\pdfoutput=1
\documentclass[10pt,letterpaper,twoside]{book}

\usepackage[utf8]{inputenc}
\usepackage[T1]{fontenc}
\usepackage[provide=*,english]{babel}
\IfFileExists{kpfonts.sty}{\usepackage{kpfonts}}{}
\IfFileExists{bbding.sty}{\usepackage{bbding}}{}
\usepackage{textcomp}
\usepackage{gensymb}

\usepackage{amsthm}
\usepackage{amsmath}
\allowdisplaybreaks
\usepackage{amssymb}
\usepackage{mathtools}
\usepackage{mathrsfs}
\IfFileExists{mathdots.sty}{\usepackage{mathdots}}{}
\IfFileExists{yhmath.sty}{\usepackage[minimal,nowidering]{yhmath}}{}
\makeatletter
\@ifpackageloaded{yhmath}{%
  \protected\edef\book@native@widetilde#1{\mathaccent"0\hexnumber@\symyhlargesymbols65\relax{#1}}%
  \protected\edef\widehat#1{\mathaccent"0\hexnumber@\symyhlargesymbols62\relax{#1}}%
  \protected\edef\widetriangle#1{\mathaccent"0\hexnumber@\symyhlargesymbols E6\relax{#1}}%
  \protected\edef\wideparen#1{\mathaccent"0\hexnumber@\symyhlargesymbols F3\relax{#1}}%
  \mathchardef\book@widetilde@glyph="0\hexnumber@\symyhlargesymbols E5\relax
  \protected\def\widetilde#1{\mathpalette\book@widetilde{#1}}%
  \def\book@widetilde#1#2{%
    \begingroup
    \setbox\z@=\hbox{$\m@th#1#2$}%
    \setbox\tw@=\hbox{$\m@th#1\book@widetilde@glyph$}%
    \ifdim\wd\z@>\wd\tw@
      \setbox4=\hbox{$\m@th#1x$}%
      \dimen@=\ht4
      \ifdim\ht\z@<\dimen@ \dimen@=\ht\z@\fi
      \mathord{\vbox{\offinterlineskip
        \hbox{\resizebox{\wd\z@}{\ht\tw@}{\copy\tw@}}%
        \kern-\dimen@
        \box\z@
      }}%
    \else
      \book@native@widetilde{#2}%
    \fi
    \endgroup
  }%
}{}
\makeatother
\usepackage{cancel}
\IfFileExists{extarrows.sty}{\usepackage{extarrows}}{}

\usepackage{xcolor}
\usepackage{graphicx}
\usepackage{pdfpages}
\usepackage{tikz}
\usepackage{tikz-cd}
\usepackage{float}
\usepackage{wrapfig}

\usepackage{tabularx}
\usepackage{longtable}
\usepackage{booktabs}
\usepackage{multirow}
\usepackage{array}

\usepackage{titlesec}
\usepackage{bold-extra}
\usepackage{ragged2e}
\usepackage{setspace}
\usepackage{enumitem}
\usepackage{tcolorbox}
\usepackage{fancyhdr}
\usepackage{geometry}

\usepackage{epigraph}
\usepackage[toc,page]{appendix}

\usepackage{caption}
\usepackage{subcaption}
\usepackage{etoolbox}

\usepackage[nocfg]{nomencl}
\IfFileExists{tracklang.sty}{%
  \usepackage[acronym,toc]{glossaries}%
}{%
  \providecommand{\makeglossaries}{}%
  \newcommand{\glsfallbackentries}{}%
  \csname define@key\endcsname{glsfallback}{name}{\def\glsfallbackname{##1}}%
  \csname define@key\endcsname{glsfallback}{description}{\def\glsfallbackdescription{##1}}%
  \csname define@key\endcsname{glsfallback}{sort}{}%
  \csname define@key\endcsname{glsfallback}{see}{}%
  \newcommand{\glsfallbackrender}[1]{%
    \def\glsfallbackname{}%
    \def\glsfallbackdescription{}%
    \setkeys{glsfallback}{##1}%
    \par\medskip\noindent
    \textbf{\glsfallbackname.}\enspace\glsfallbackdescription\par
  }%
  \newcommand{\newglossaryentry}[2]{%
    \gappto\glsfallbackentries{\glsfallbackrender{##2}}%
  }%
  \providecommand{\glsadd}[1]{}%
  \newcommand{\printglossary}[1][]{%
    \cleardoublepage
    \chapter*{Glossary}%
    \markboth{Glossary}{Glossary}%
    \addcontentsline{toc}{chapter}{Glossary}%
    \glsfallbackentries
  }%
}
\usepackage{imakeidx}

\usepackage[
  colorlinks=false,
  pdfborder={0 0 1},
  linkbordercolor={0 0 1},
  pdftitle={Global analysis. An introduction to nonlinear analysis and its variational methods on Riemannian manifolds},
  pdfauthor={Carlos Daniel Velázquez Mendoza; María de los Ángeles Sandoval Romero; Romulo Diaz Carlos}
]{hyperref}
\IfFileExists{bookmark.sty}{\usepackage{bookmark}}{}
\usepackage{cleveref}

\graphicspath{{figuras/}{figuras/imagenes/}}
\usetikzlibrary{babel}
\tcbuselibrary{skins,breakable,hooks}
\definecolor{miColorPrincipal}{RGB}{128,0,32}

\newenvironment{figura}[1][htbp]{%
    \begin{figure}[#1]
    \centering
    \begin{tcolorbox}[
        enhanced,
        colback=miColorPrincipal!4,
        frame hidden,
        boxrule=0pt,
        borderline west={0.5mm}{0mm}{miColorPrincipal},
        left=4mm,
        right=4mm,
        top=3mm,
        bottom=2mm,
        arc=.5mm
    ]
    \centering
}{%
    \end{tcolorbox}
    \end{figure}
}
\newtcolorbox{notaFlotanteBox}[1]{%
  enhanced,
  colback=miColorPrincipal!6,
  frame hidden,
  boxrule=0pt,
  borderline west={0.8pt}{0pt}{miColorPrincipal!85!black},
  left=8pt,right=6pt,top=6pt,bottom=6pt,
  arc=2pt,
  title={#1},
  coltitle=miColorPrincipal!85!black,
  fonttitle=\normalfont\small\sffamily,
  colbacktitle=miColorPrincipal!0,
  titlerule=0pt,
  toptitle=2pt,
  bottomtitle=4pt,
  fontupper=\normalfont\small,
}

\newcommand{\epigrafeFancy}[2]{%
    \vspace{0.5cm} 
    \begin{flushright} 
        \begin{tcolorbox}[
            blanker,       
            width=0.7\textwidth, 
            borderline west={3pt}{0pt}{miColorPrincipal}, 
            left=12pt,     
            top=2pt, bottom=2pt,
            fontupper=\small\itshape\color{black!85}, 
        ]
            ``#1''

            \par\vspace{8pt} 
            \hfill \normalfont\scshape\bfseries --- #2 
        \end{tcolorbox}
    \end{flushright}
    \vspace{1.5cm} 
}
\newcommand{\cajaTitulo}[1]{%
  \begingroup
    \setlength{\fboxsep}{3pt}%
    \colorbox{miColorPrincipal}{%
      \hspace{1pt}\textcolor{white}{\bfseries\sffamily\Large #1}\hspace{1pt}}%
  \endgroup
}
\titleformat{\part}[display]
  {\thispagestyle{empty}\raggedleft}
  {
    \tcbox[
        on line,
        sharp corners=all,
        enhanced,
        colback=miColorPrincipal,
        colframe=miColorPrincipal,
        colupper=white,
        boxrule=0pt,
        boxsep=5pt,
        left=4pt,
        right=4pt,
        top=3pt,
        bottom=3pt,
        fontupper=\bfseries\sffamily\fontsize{30}{30}\selectfont
    ]{\partname\ \thepart}%
  }
  {20pt}
  {\sffamily\bfseries\fontsize{50}{60}\selectfont\color{black}}
  [\vspace{20pt}{\color{miColorPrincipal}\titlerule[3pt]}]

\titleformat{\chapter}[display]
  {\sffamily\bfseries\color{black}} 
  {\flushright\fontsize{50}{50}\selectfont\color{miColorPrincipal}\thechapter}
  {-10pt}
  {\flushright\Huge}
  [\vspace{5pt}{\color{miColorPrincipal}\titlerule[2pt]}]
\titleformat{\section}
  {\Large\sffamily\bfseries\color{miColorPrincipal}} 
  {\cajaTitulo{\thesection}} 
  {10pt} 
  {}

\titleformat{\subsection}
  {\large\sffamily\bfseries\color{black}}
  {\color{miColorPrincipal}\thesubsection}
  {1em}
  {}

\fancyhfoffset{0pt}

\newcommand{\cajaPagina}[1]{%
  \begingroup
    \setlength{\fboxsep}{2pt}%
    \colorbox{miColorPrincipal}{%
      \hspace{2pt}\textcolor{white}{\bfseries\sffamily\small #1}\hspace{2pt}}%
  \endgroup
}

\renewcommand{\headrulewidth}{1pt}
\renewcommand{\headrule}{\hbox to\headwidth{%
    \color{miColorPrincipal}\leaders\hrule height \headrulewidth\hfill}}

\newcommand{\estiloHeader}[1]{{\sffamily\color{black}\small #1}}

\makeindex

\newtcolorbox{semblanzaHistorica}[1]{%
  enhanced,
  breakable,
  colback=miColorPrincipal!3,
  colframe=miColorPrincipal!55!black,
  boxrule=0.35pt,
  borderline west={1.1pt}{0pt}{miColorPrincipal},
  left=8pt,
  right=8pt,
  top=7pt,
  bottom=7pt,
  arc=1.5pt,
  title={#1},
  coltitle=miColorPrincipal!90!black,
  fonttitle=\sffamily\bfseries,
  colbacktitle=miColorPrincipal!3,
  titlerule=0pt,
  before skip=10pt,
  after skip=10pt
}

\newtheoremstyle{theorem}
    {0pt}{0pt}{\normalfont}{0pt}
    {}{}{0.25em}
    {{\sffamily\bfseries\color{theoremcolor}\thmname{#1}~\thmnumber{\textup{#2}}.}
        \thmnote{\normalfont\color{black}~(#3)}}

\newtheoremstyle{definition}
    {0pt}{0pt}{\normalfont}{0pt}
    {}{}{0.25em}
    {{\sffamily\bfseries\color{definitioncolor}\thmname{#1}~\thmnumber{\textup{#2}}.}
        \thmnote{\normalfont\color{black}~(#3)}}

\newtheoremstyle{proposition}
    {0pt}{0pt}{\normalfont}{0pt}
    {}{}{0.25em}
    {{\sffamily\bfseries\color{propositioncolor}\thmname{#1}~\thmnumber{\textup{#2}}.}
        \thmnote{\normalfont\color{black}~(#3)}}

\newtheoremstyle{corollary}
    {0pt}{0pt}{\normalfont}{0pt}
    {}{}{0.25em}
    {{\sffamily\bfseries\color{corollarycolor}\thmname{#1}~\thmnumber{\textup{#2}}.}
        \thmnote{\normalfont\color{black}~(#3)}}

\newtheoremstyle{lemma}
    {0pt}{0pt}{\normalfont}{0pt}
    {}{}{0.25em}
    {{\sffamily\bfseries\color{lemmacolor}\thmname{#1}~\thmnumber{\textup{#2}}.}
        \thmnote{\normalfont\color{black}~(#3)}}

\newtheoremstyle{remark}
    {0pt}{0pt}{\normalfont}{0pt}
    {}{}{0.25em}
    {{\sffamily\bfseries\color{remarkcolor}\thmname{#1}~\thmnumber{\textup{#2}}.}
        \thmnote{\normalfont\color{black}~(#3)}}

\newtheoremstyle{note}
    {0pt}{0pt}{\normalfont}{0pt}
    {}{}{0.25em}
    {{\sffamily\bfseries\color{notecolor}\thmname{#1}~\thmnumber{\textup{#2}}.}
        \thmnote{\normalfont\color{black}~(#3)}}

\newtheoremstyle{example}
    {0pt}{0pt}{\normalfont}{0pt}
    {}{}{0.25em}
    {{\sffamily\bfseries\color{examplecolor}\thmname{#1}~\thmnumber{\textup{#2}}.}
        \thmnote{\normalfont\color{black}~(#3)}}
\newtheoremstyle{notation}
    {0pt}{0pt}{\normalfont}{0pt}
    {}{}{0.25em}
    {{\sffamily\bfseries\color{notationcolor}\thmname{#1}~\thmnumber{\textup{#2}}.}
        \thmnote{\normalfont\color{black}~(#3)}}
\theoremstyle{theorem}
\newtheorem{theorem}{Theorem}[section]

\theoremstyle{definition}
\newtheorem{definition}[theorem]{Definition}

\theoremstyle{proposition}
\newtheorem{proposition}[theorem]{Proposition}

\theoremstyle{corollary}
\newtheorem{corollary}[theorem]{Corollary}

\theoremstyle{lemma}
\newtheorem{lemma}[theorem]{Lemma}

\theoremstyle{remark}
\newtheorem{remark}[theorem]{Remark}

\theoremstyle{note}
\newtheorem{note}[theorem]{Note}
\theoremstyle{notation}
\newtheorem{notation}[theorem]{Notation}
\theoremstyle{example}
\newtheorem{example}[theorem]{Example}

\colorlet{exercisecolor}{green!50!blue}

\newtheoremstyle{exercise}%
    {0pt}{0pt}
    {\normalfont}
    {0pt}
    {}{}
    {0.25em}
    {{\sffamily\bfseries\color{exercisecolor}\thmname{#1}~\thmnumber{\textup{#2}}.}%
     \thmnote{\normalfont\color{black}~(#3)}}

\theoremstyle{exercise}
\newtheorem{exercise}[theorem]{Exercise}

\tcolorboxenvironment{exercise}{
    enhanced, pad at break*=1mm, breakable,
    left=4mm, right=4mm, top=1mm, bottom=1mm,
    colback=exercisecolor!10, boxrule=0pt, frame hidden,
    borderline west={0.5mm}{0mm}{exercisecolor}, arc=.5mm
}

\colorlet{theoremcolor}{brown!50!blue}
\colorlet{definitioncolor}{purple!50!brown}
\colorlet{propositioncolor}{blue!60!brown}
\colorlet{corollarycolor}{purple!60!blue}
\colorlet{lemmacolor}{brown!60!purple}
\colorlet{remarkcolor}{blue!40!brown}
\colorlet{notecolor}{brown!50!purple}
\colorlet{examplecolor}{purple!50!blue}
\colorlet{notationcolor}{gray!60!black}
\tcolorboxenvironment{theorem}{
    enhanced, pad at break*=1mm, breakable,
    left=4mm, right=4mm, top=1mm, bottom=1mm,
    colback=theoremcolor!10, boxrule=0pt, frame hidden,
    borderline west={0.5mm}{0mm}{theoremcolor}, arc=.5mm
}

\tcolorboxenvironment{definition}{
    enhanced, pad at break*=1mm, breakable,
    left=4mm, right=4mm, top=1mm, bottom=1mm,
    colback=definitioncolor!10, boxrule=0pt, frame hidden,
    borderline west={0.5mm}{0mm}{definitioncolor}, arc=.5mm
}
\tcolorboxenvironment{notation}{
    enhanced, pad at break*=1mm, breakable,
    left=4mm, right=4mm, top=1mm, bottom=1mm,
    colback=notationcolor!8,
    boxrule=0pt, frame hidden,
    borderline west={0.5mm}{0mm}{notationcolor},
    arc=.5mm
}
\tcolorboxenvironment{proposition}{
    enhanced, pad at break*=1mm, breakable,
    left=4mm, right=4mm, top=1mm, bottom=1mm,
    colback=propositioncolor!10, boxrule=0pt, frame hidden,
    borderline west={0.5mm}{0mm}{propositioncolor}, arc=.5mm
}

\tcolorboxenvironment{corollary}{
    enhanced, pad at break*=1mm, breakable,
    left=4mm, right=4mm, top=1mm, bottom=1mm,
    colback=corollarycolor!10, boxrule=0pt, frame hidden,
    borderline west={0.5mm}{0mm}{corollarycolor}, arc=.5mm
}

\tcolorboxenvironment{lemma}{
    enhanced, pad at break*=1mm, breakable,
    left=4mm, right=4mm, top=1mm, bottom=1mm,
    colback=lemmacolor!10, boxrule=0pt, frame hidden,
    borderline west={0.5mm}{0mm}{lemmacolor}, arc=.5mm
}

\tcolorboxenvironment{remark}{
    enhanced, pad at break*=1mm, breakable,
    left=4mm, right=4mm, top=1mm, bottom=1mm,
    colback=remarkcolor!10, boxrule=0pt, frame hidden,
    borderline west={0.5mm}{0mm}{remarkcolor}, arc=.5mm
}

\tcolorboxenvironment{note}{
    enhanced, pad at break*=1mm, breakable,
    left=4mm, right=4mm, top=1mm, bottom=1mm,
    colback=notecolor!10, boxrule=0pt, frame hidden,
    borderline west={0.5mm}{0mm}{notecolor}, arc=.5mm
}

\tcolorboxenvironment{example}{
    enhanced, pad at break*=1mm, breakable,
    left=4mm, right=4mm, top=1mm, bottom=1mm,
    colback=white, 
    boxrule=0pt, frame hidden,
    borderline west={0.5mm}{0mm}{examplecolor}, arc=.5mm
}

\AddToHook{env/theorem/after}{\colorlet{proofcolor}{theoremcolor}}
\AddToHook{env/definition/after}{\colorlet{proofcolor}{definitioncolor}}
\AddToHook{env/proposition/after}{\colorlet{proofcolor}{propositioncolor}}
\AddToHook{env/corollary/after}{\colorlet{proofcolor}{corollarycolor}}
\AddToHook{env/lemma/after}{\colorlet{proofcolor}{lemmacolor}}
\AddToHook{env/remark/after}{\colorlet{proofcolor}{remarkcolor}}
\AddToHook{env/note/after}{\colorlet{proofcolor}{notecolor}}
\AddToHook{env/example/after}{\colorlet{proofcolor}{examplecolor}}
\AddToHook{env/notation/after}{\colorlet{proofcolor}{notationcolor}}

\let\qedsymbolMyOriginal\qedsymbol
\renewcommand{\qedsymbol}{%
    \color{proofcolor}\qedsymbolMyOriginal%
}
\newtheoremstyle{proofstyle}
    {0pt}{0pt}{\normalfont}{0pt}
    {}{}{0.25em}
    {{\sffamily\bfseries\color{proofcolor}\thmname{#1}~\thmnumber{\textup{#2}}.}
        \thmnote{\normalfont\color{black}~(#3)}}

\theoremstyle{proofstyle}

\colorlet{proofcolor}{blue!60!brown}

\tcolorboxenvironment{proof}{
    enhanced, pad at break*=1mm, breakable,
    left=4mm, right=4mm, top=1mm, bottom=1mm,
    colback=white,  
    boxrule=0pt, frame hidden,
    borderline west={0.5mm}{0mm}{proofcolor},
    arc=.5mm
}

\AddToHook{env/proof/after}{\colorlet{proofcolor}{proofcolor}}

\let\qedsymbolProofOriginal\qedsymbol
\renewcommand{\qedsymbol}{%
    \color{proofcolor}\qedsymbolProofOriginal%
}

\makeatletter
\def\smallunderbrace#1{\mathop{\vtop{\m@th\ialign{##\crcr
   $\hfil\displaystyle{#1}\hfil$\crcr
   \noalign{\kern3\p@\nointerlineskip}%
   \tiny\upbracefill\crcr\noalign{\kern3\p@}}}}\limits}
\makeatother
\makeatletter
\newenvironment{solution}[1][\relax]{%
  \def\solutionname{%
    \normalfont\sffamily\bfseries%
    \color{proofcolor}Solution%
  }%
  \ifx\relax#1\relax
    \proof[\solutionname]%
  \else
    \proof[\solutionname\ (#1)]%
  \fi
}{%
  \endproof
}
\makeatother

\tcolorboxenvironment{solution}{
    enhanced, pad at break*=0mm, breakable,
    left=4mm, right=4mm, top=1mm, bottom=1mm,
    colback=white,
    boxrule=0pt, frame hidden,
    borderline west={0.5mm}{0mm}{proofcolor},
    arc=.5mm
}

\DeclareMathOperator{\sign}{sign}
\newcommand\supp{\mathrm{supp}}

\renewcommand{\phi}{\varphi}

\addto\captionsenglish{}

\makenomenclature

\renewcommand{\nomgroup}[1]{%
  \item[\bfseries
    \if A#1 Sets and general conventions%
    \else\if B#1 Topology and functional analysis%
    \else\if C#1 Differential geometry%
    \else\if D#1 Riemannian geometry%
    \else\if E#1 Vector bundles and sections%
    \else\if F#1 Sobolev spaces%
    \else\if G#1 Differential operators%
    \else\if H#1 Infinite-dimensional geometry%
    \else\if I#1 Variational methods%
    \else #1%
    \fi\fi\fi\fi\fi\fi\fi\fi\fi]}

\let\emptyset\varnothing
\newcommand\underrel[3][]{\mathrel{\mathop{#3}\limits_{%
      \ifx c#1\relax\mathclap{#2}\else#2\fi}}}

\makeglossaries

\newglossaryentry{variedad-riemanniana}{
  name={Riemannian manifold},
  sort={Riemannian manifold},
  description={A smooth manifold equipped with a Riemannian metric; see Definition~\ref{def:variedades-riemannianas-variedad}},
  see={metrica-riemanniana,medida-riemann-lebesgue}
}

\newglossaryentry{metrica-riemanniana}{
  name={Riemannian metric},
  sort={Riemannian metric},
  description={A metric tensor field that determines the inner product on each tangent space; see Definition~\ref{def:variedades-riemannianas-variedad}},
  see={variedad-riemanniana}
}

\newglossaryentry{medida-riemann-lebesgue}{
  name={Riemann--Lebesgue measure},
  sort={Riemann Lebesgue measure},
  description={The canonical measure associated with a Riemannian metric, used to formulate global integrals without assuming orientability; see Definition~\ref{medibles ajenos dos a dos, definicion de medida Riemann--Lebesgue}},
  see={variedad-riemanniana}
}

\newglossaryentry{haz-vectorial}{
  name={smooth vector bundle},
  sort={smooth vector bundle},
  description={A vector bundle whose projection, total space, and local trivializations have a smooth structure; see Definition~\ref{def:nociones-fundamentales-haz-vectorial-suave-de-rango-sobre}},
  see={metrica-fibrada,conexion}
}

\newglossaryentry{metrica-fibrada}{
  name={bundle metric},
  sort={bundle metric},
  description={A smooth choice of an inner product on each fiber of a vector bundle; see Definition~\ref{def producto fibrado}},
  see={haz-vectorial,conexion-compatible}
}

\newglossaryentry{conexion}{
  name={connection},
  sort={connection},
  description={An operator that differentiates sections of a bundle in tangent directions and satisfies linearity and the Leibniz rule; see Definition~\ref{definicion 1 conexion}},
  see={derivada-covariante,conexion-compatible}
}

\newglossaryentry{conexion-compatible}{
  name={compatible connection},
  sort={compatible connection},
  description={A connection that preserves the bundle metric under differentiation; see Definition~\ref{conexion compatible con metrica}},
  see={conexion,metrica-fibrada}
}

\newglossaryentry{derivada-covariante}{
  name={covariant derivative},
  sort={covariant derivative},
  description={The derivative of a section determined by a connection and a tangent direction},
  see={conexion,derivada-covariante-debil}
}

\newglossaryentry{derivada-covariante-debil}{
  name={weak covariant derivative},
  sort={weak covariant derivative},
  description={The distributional extension of the covariant derivative, defined through integration by parts; see Definition~\ref{def:definicion-sobolev-haces-derivada-variedad-riemanniana-haz-vectorial}},
  see={derivada-covariante,espacio-sobolev}
}

\newglossaryentry{adjunto-formal}{
  name={formal adjoint operator},
  sort={formal adjoint operator},
  description={A differential operator characterized by the integration-by-parts identity with respect to the metrics and the Riemann--Lebesgue measure; see Definition~\ref{def: adjunto formal de un pdo}},
  see={laplaciano-bochner,laplaciano-hodge}
}

\newglossaryentry{operador-traza}{
  name={trace operator},
  sort={trace operator},
  description={The continuous extension of restriction to the boundary for Sobolev classes of sufficient regularity; see Theorem~\ref{teo:traza-entera-haces}},
  see={espacio-sobolev}
}

\newglossaryentry{espacio-sobolev}{
  name={Sobolev space},
  sort={Sobolev space},
  description={A function space that controls a function or section and its weak derivatives up to a given order; see Definitions~\ref{def:sobolev-en-abiertos-euclidianos-espacio-de-sobolev} and~\ref{def:definicion-sobolev-haces-norma-variedad-riemanniana-haz-vectorial}},
  see={derivada-covariante-debil,operador-traza}
}

\newglossaryentry{desigualdad-kato}{
  name={Kato inequality},
  sort={Kato inequality},
  description={An estimate comparing the weak derivative of the norm of a section with its covariant derivative; see Theorem~\ref{desigualdad de kato}},
  see={derivada-covariante-debil}
}

\newglossaryentry{desigualdad-poincare}{
  name={Poincaré inequality},
  sort={Poincare inequality},
  description={An estimate for a function or section in terms of its derivatives, under the corresponding geometric or boundary conditions; see Theorems~\ref{desigualdad de poincaré} and~\ref{teo:poincare-modulo-nucleo-conexion}},
  see={espacio-sobolev}
}

\newglossaryentry{geometria-acotada}{
  name={bounded geometry},
  sort={bounded geometry},
  description={A uniform geometric condition formulated in terms of the injectivity radius and bounds on the derivatives of curvature; see Definition~\ref{def: geometria acotada}},
  see={variedad-riemanniana}
}

\newglossaryentry{laplaciano-bochner}{
  name={Bochner Laplacian},
  sort={Bochner Laplacian},
  description={A second-order operator constructed from a connection and its formal adjoint; see Definition~\ref{def:operadores-diferenciales-en-haces-laplaciano-de-bochner}},
  see={adjunto-formal,laplaciano-hodge}
}

\newglossaryentry{laplaciano-hodge}{
  name={Hodge--de Rham Laplacian},
  sort={Hodge de Rham Laplacian},
  description={An operator on differential forms given by the sum of the compositions of the exterior derivative and its formal adjoint; see Definition~\ref{def:operadores-diferenciales-en-haces-laplaciano-de-hodge-derham}},
  see={adjunto-formal,laplaciano-bochner}
}

\newglossaryentry{distribucion-haz}{
  name={distribution in a vector bundle},
  sort={distribution in a vector bundle},
  description={A continuous linear operator defined on test sections of the dual bundle; see Definition~\ref{def:distribuciones-en-haces-distribuciones}},
  see={haz-vectorial}
}

\newglossaryentry{derivada-frechet}{
  name={Fréchet derivative},
  sort={Frechet derivative},
  description={The first-order continuous linear approximation to a function between Banach spaces; see Definition~\ref{def:derivada-frechet}},
  see={variedad-banach}
}

\newglossaryentry{variedad-banach}{
  name={Banach manifold},
  sort={Banach manifold},
  description={A manifold whose charts take values in open subsets of Banach spaces and whose coordinate changes are Fréchet differentiable; see Definition~\ref{def:variedad-banach}},
  see={derivada-frechet,subvariedad-banach}
}

\newglossaryentry{subvariedad-banach}{
  name={Banach submanifold},
  sort={Banach submanifold},
  description={A subset that can be locally straightened onto a closed complemented subspace of the model space; see Definition~\ref{def:subvariedad-encajada-banach}},
  see={variedad-banach,nivel-regular}
}

\newglossaryentry{inmersion-banach}{
  name={immersion between Banach manifolds},
  sort={immersion between Banach manifolds},
  description={A map whose differential admits a continuous linear left inverse; see Definition~\ref{def:inmersion-submersion-split-banach}},
  see={variedad-banach,subvariedad-banach}
}

\newglossaryentry{submersion-banach}{
  name={submersion between Banach manifolds},
  sort={submersion between Banach manifolds},
  description={A map whose differential admits a continuous linear right inverse; see Definition~\ref{def:inmersion-submersion-split-banach}},
  see={variedad-banach,nivel-regular}
}

\newglossaryentry{valor-regular-banach}{
  name={split regular value},
  sort={split regular value},
  description={A value whose fiber consists of points at which the differential admits a continuous linear right inverse; see Definition~\ref{def:valor-regular-split-banach}},
  see={submersion-banach,subvariedad-banach}
}

\newglossaryentry{funcion-flan-banach}{
  name={bump function on a Banach space},
  sort={bump function on a Banach space},
  description={A nonzero differentiable function with bounded support; see Definition~\ref{def:funcion-flan-banach}},
  see={variedad-banach,particion-unidad-banach}
}

\newglossaryentry{particion-unidad-banach}{
  name={smooth partition of unity on a Banach manifold},
  sort={smooth partition of unity on a Banach manifold},
  description={A locally finite family of nonnegative smooth functions subordinate to an open cover and with sum equal to one; see Theorem~\ref{teo:particiones-unidad-variedad-banach}},
  see={funcion-flan-banach,variedad-banach}
}

\newglossaryentry{nivel-regular}{
  name={regular level set},
  sort={regular level set},
  description={A level set of a functional whose derivative does not vanish at its points and which, in the setting developed here, has a submanifold structure; see Theorem~\ref{teo:conjuntos-nivel-regulares-banach}},
  see={subvariedad-banach,multiplicador-lagrange}
}

\newglossaryentry{funcion-caratheodory}{
  name={Carathéodory function},
  sort={Caratheodory function},
  description={A function that is measurable in the geometric variable and continuous in the real variable; see Definition~\ref{def:funcion-caratheodory-variedad}},
  see={solucion-debil,funcional-energia}
}

\newglossaryentry{solucion-debil}{
  name={weak solution},
  sort={weak solution},
  description={An element of the function space that satisfies the integral identity obtained by testing the equation against admissible functions; see Definition~\ref{def:solucion-debil-problema-modelo}},
  see={funcional-energia,metodo-variacional}
}

\newglossaryentry{funcional-energia}{
  name={energy functional},
  sort={energy functional},
  description={A variational functional whose critical points represent weak solutions of the associated Euler--Lagrange equation; see Section~\ref{sec:funcional-energia-puntos-criticos}},
  see={punto-critico,metodo-variacional}
}

\newglossaryentry{punto-critico}{
  name={critical point},
  sort={critical point},
  description={A point at which the Fréchet derivative of the functional under consideration vanishes},
  see={derivada-frechet,funcional-energia}
}

\newglossaryentry{metodo-variacional}{
  name={variational method},
  sort={variational method},
  description={A procedure that translates a differential equation into a critical point problem for a functional; see Section~\ref{sec:funcional-energia-puntos-criticos}},
  see={funcional-energia,metodo-directo}
}

\newglossaryentry{semicontinuidad-inferior}{
  name={lower semicontinuity},
  sort={lower semicontinuity},
  description={The property of a functional expressing stability from below with respect to the topology under consideration; see Definitions~\ref{def:semicontinuidad-inferior} and~\ref{def:semicontinuidad-inferior-debil}},
  see={metodo-directo}
}

\newglossaryentry{coercividad}{
  name={coercivity},
  sort={coercivity},
  description={The property of a functional whose value tends to \(+\infty\) as the norm of its argument tends to \(+\infty\)},
  see={metodo-directo}
}

\newglossaryentry{metodo-directo}{
  name={direct method in the calculus of variations},
  sort={direct method in the calculus of variations},
  description={An existence principle for minimizers based on weak compactness, coercivity, and lower semicontinuity; see Theorem~\ref{teo:metodo-directo-reflexivo}},
  see={semicontinuidad-inferior,coercividad}
}

\newglossaryentry{multiplicador-lagrange}{
  name={Lagrange multiplier},
  sort={Lagrange multiplier},
  description={A scalar expressing the dependence between the derivative of a functional and the derivative of a constraint at a constrained critical point; see Corollary~\ref{teo:multiplicadores-lagrange-banach}},
  see={nivel-regular,variedad-nehari}
}

\newglossaryentry{variedad-nehari}{
  name={Nehari manifold},
  sort={Nehari manifold},
  description={The set of nonzero elements defined by the condition \(DI(u)(u)=0\), used as a natural constraint for certain functionals; see Section~\ref{sec:metodo-nehari}},
  see={multiplicador-lagrange,funcional-energia}
}

\newglossaryentry{distribucion-temperada}{
  name={tempered distribution},
  sort={tempered distribution},
  description={A continuous linear functional on the Schwartz class; see Definition~\ref{def:distribucion-temperada}},
  see={transformada-fourier}
}
\newglossaryentry{transformada-fourier}{
  name={Fourier transform},
  sort={Fourier transform},
  description={An operator that decomposes a function or tempered distribution into frequencies; see Definitions~\ref{def:transformada-fourier-L1} and~\ref{def:fourier-distribucion-temperada}},
  see={distribucion-temperada}
}
\newglossaryentry{interpolacion-real}{
  name={real interpolation}, sort={real interpolation},
  description={A method for constructing intermediate spaces using the $K$- and $J$-functionals; see Definition~\ref{def:espacio-interpolacion-real-K}},
  see={interpolacion-compleja}
}
\newglossaryentry{interpolacion-compleja}{
  name={complex interpolation}, sort={complex interpolation},
  description={Calderón's method using holomorphic functions on a strip; see Definition~\ref{def:espacio-calderon-interpolacion-compleja}},
  see={interpolacion-real}
}
\newglossaryentry{espacio-slobodeckij}{
  name={Sobolev--Slobodeckij space}, sort={Sobolev Slobodeckij space},
  description={A space of noninteger order defined using differences and the Gagliardo seminorm; see Definition~\ref{def:seminorma-gagliardo-slobodeckij}},
  see={espacio-bessel}
}
\newglossaryentry{espacio-bessel}{
  name={Bessel potential space}, sort={Bessel potential space},
  description={The space $H^{s,p}$ defined by the multiplier $(1+\|\xi\|^2)^{\frac{s}{2}}$; see Definition~\ref{def:regularidad-intermedia-norma-espacio-sobolev-fraccionario}},
  see={espacio-slobodeckij}
}

\newglossaryentry{variedad-hilbert}{
  name={Hilbert manifold},
  sort={Hilbert manifold},
  description={A Banach manifold modeled on a Hilbert space; see Definition~\ref{def:variedad-hilbert}},
  see={variedad-banach,metrica-riemanniana-fuerte}
}

\newglossaryentry{metrica-riemanniana-fuerte}{
  name={strong Riemannian metric},
  sort={strong Riemannian metric},
  description={A smooth metric whose tangent norm induces the Hilbert topology and whose musical morphism is an isomorphism; see Section~\ref{sec:metricas-hilbert-fuertes-debiles}},
  see={variedad-hilbert,banach-finsler}
}

\newglossaryentry{banach-finsler}{
  name={Banach--Finsler structure},
  sort={Banach Finsler structure},
  description={A continuous family of complete tangent norms satisfying the Palais condition: local $(1+\varepsilon)$-comparability in charts for every $\varepsilon>0$; see Definition~\ref{def:estructura-banach-finsler}},
  see={metrica-riemanniana-fuerte,palais-smale}
}

\newglossaryentry{palais-smale}{
  name={Palais--Smale condition},
  sort={Palais Smale condition},
  description={A compactness condition for sequences with bounded energy whose derivative tends to zero; see Definitions~\ref{def:palais-smale-banach} and~\ref{def:palais-smale-finsler}},
  see={pseudogradiente,banach-finsler}
}

\newglossaryentry{pseudogradiente}{
  name={pseudogradient vector field},
  sort={pseudogradient vector field},
  description={A locally Lipschitz vector field controlled by the norm of the derivative and producing uniform descent of the functional; see Definition~\ref{def:pseudogradiente-banach}},
  see={palais-smale}
}

\newglossaryentry{uniformemente-regular}{
  name={uniformly regular manifold},
  sort={uniformly regular manifold},
  description={A manifold equipped with a class of normalized, shrinkable atlases of finite multiplicity, with uniformly controlled coordinate changes in the sense of Amann; see Definition~\ref{def:estructura-uniformemente-regular-amann}},
  see={geometria-acotada,estructura-uniforme}
}

\newglossaryentry{estructura-uniforme}{
  name={uniform structure},
  sort={uniform structure},
  description={A filter of entourages of the diagonal that makes it possible to define uniform continuity, Cauchy filters, and completeness; see Definition~\ref{def:estructura-uniforme}},
  see={uniformemente-regular}
}

\newglossaryentry{mapeos-sobolev}{
  name={manifold of Sobolev maps},
  sort={manifold of Sobolev maps},
  description={A component of Sobolev maps modeled on sections of $f^*TN$ through exponential charts; see Theorem~\ref{teo:cartas-exponenciales-mapeos-sobolev}},
  see={variedad-banach,escala-hilbert}
}

\newglossaryentry{escala-hilbert}{
  name={Hilbert scale},
  sort={Hilbert scale},
  description={A chain of Hilbert spaces with continuous dense inclusions that records Sobolev regularity; see Definition~\ref{def:escala-hilbert}},
  see={variedad-ilh,mapeos-sobolev}
}

\newglossaryentry{variedad-ilh}{
  name={ILH manifold},
  sort={ILH manifold},
  description={The smooth inverse limit of a compatible tower of Hilbert manifolds; see Definition~\ref{def:variedad-grupo-ilh}},
  see={escala-hilbert}
}

\newglossaryentry{espacios-bc-buc}{
  name={spaces $BC^j$ and $BUC^j$},
  sort={spaces BC^j and BUC^j},
  description={Spaces of $C^j$ functions whose derivatives are bounded and, in the second case, uniformly continuous; see Definition~\ref{def:espacios-bc-buc}},
  see={uniformemente-regular}
}

\newglossaryentry{dual-topologico-banach}{
  name={continuous dual of a Banach space},
  sort={continuous dual of a Banach space},
  description={The space of continuous linear functionals, denoted by $X'$; see Proposition~\ref{prop:dual-topologico-banach}},
  see={haz-dual-topologico-banach,cotangente-banach}
}

\newglossaryentry{haz-dual-topologico-banach}{
  name={continuous dual Banach bundle},
  sort={continuous dual Banach bundle},
  description={A Banach bundle whose fibers are the continuous duals of the fibers of the original bundle; see Proposition~\ref{prop:haz-dual-banach}},
  see={dual-topologico-banach,cotangente-banach}
}

\newglossaryentry{cotangente-banach}{
  name={Banach cotangent bundle},
  sort={Banach cotangent bundle},
  description={The continuous dual bundle of the tangent bundle of a Banach manifold, denoted by $T'M$; see Definition~\ref{def:haz-cotangente-banach}},
  see={dual-topologico-banach,haz-dual-topologico-banach}
}

\newglossaryentry{metrica-riemanniana-debil}{
  name={weak Riemannian metric},
  sort={weak Riemannian metric},
  description={A smooth metric whose musical map need not be surjective onto the continuous dual; see Definition~\ref{def:variedad-hilbert}},
  see={metrica-riemanniana-fuerte,variedad-hilbert}
}

\newglossaryentry{finsler-clasico}{
  name={classical Finsler structure},
  sort={classical Finsler structure},
  description={A finite-dimensional geometry defined by a Minkowski functional that is smooth and strongly convex away from the zero section; see Definition~\ref{def:finsler-clasico}},
  see={banach-finsler}
}

\newglossaryentry{metrica-extendida}{
  name={extended metric},
  sort={extended metric},
  description={A distance function that may take the value $+\infty$ and whose restriction to each finite-distance component is a metric; see Definition~\ref{def:metrica-extendida-componente-finita}},
  see={estructura-uniforme}
}

\newglossaryentry{espectro-esencial}{
  name={essential spectrum},
  sort={essential spectrum},
  description={The part of the spectrum of a self-adjoint operator that does not consist of isolated eigenvalues of finite multiplicity; see Definition~\ref{def:espectro-esencial-brecha-cap24}},
  see={geometria-acotada}
}

\newglossaryentry{espacio-configuraciones}{
  name={configuration space},
  sort={configuration space},
  description={A function space of geometric objects equipped with the topology or differentiable structure appropriate for formulating equations and symmetries},
  see={mapeos-sobolev,variedad-ilh}
}

\newglossaryentry{operador-pseudodiferencial}{
  name={pseudodifferential operator},
  sort={pseudodifferential operator},
  description={An operator defined locally by quantizing a symbol and stable under the symbolic calculus modulo smoothing operators},
  see={parametrix,frente-onda}
}

\newglossaryentry{parametrix}{
  name={parametrix},
  sort={parametrix},
  description={A two-sided inverse of an operator modulo smoothing operators},
  see={operador-pseudodiferencial}
}

\newglossaryentry{frente-onda}{
  name={wavefront set},
  sort={wavefront set},
  description={A conic subset of the cotangent bundle with the zero section removed, identifying the points and codirections in which a distribution is not smooth},
  see={operador-pseudodiferencial,distribucion-temperada}
}

\begin{document}

\pagenumbering{Roman}
\includepdf[pages=1]{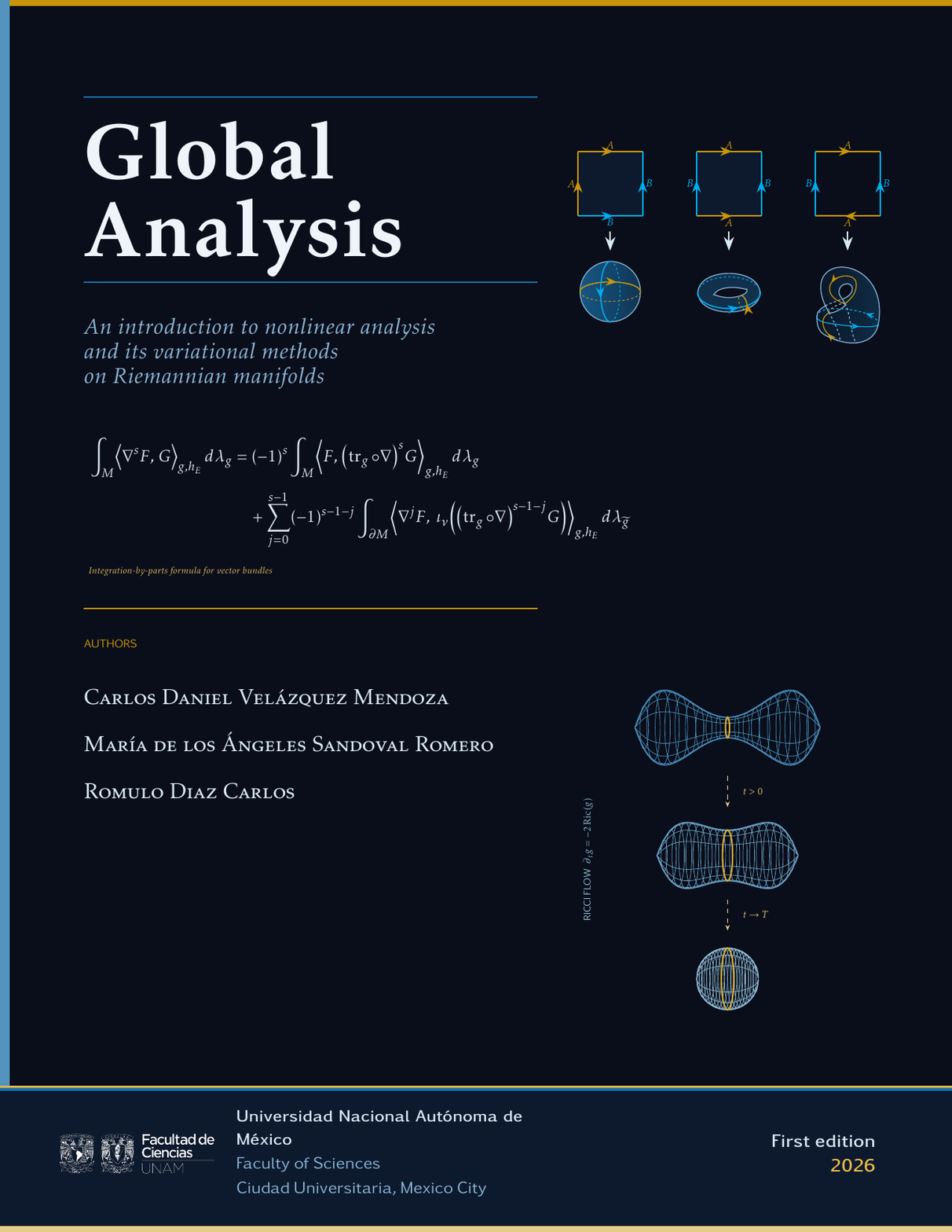}

\newpage

\setcounter{page}{2}

\

\

\

\

\

\

\epigrafeFancy{It is well known that geometry presupposes not only the concept of space but also the first fundamental notions for constructions in space as given in advance. It only gives nominal definitions for them, while the essential means of determining them appear in the form of axioms. The relationship of these presumptions is left in the dark; one sees neither whether and in how far their connection is necessary, nor a priori whether it is possible. From Euclid to Legendre, to name the most renowned of modern writers on geometry, this darkness has been lifted neither by the mathematicians nor the philosophers who have laboured upon it.}{Bernhard Riemann}





\nomenclature[A01]{$\mathbb{K}$}{Scalar field, equal to $\mathbb{R}$ or $\mathbb{C}$}
\nomenclature[A02]{$\mathbb{N}$}{Set of natural numbers}
\nomenclature[A03]{$\mathbb{N}_{0}$}{The set $\mathbb{N}\cup\{0\}$}
\nomenclature[A04]{$\mathbb{Z}$}{Set of integers}
\nomenclature[A05]{$\mathbb{R}^{n}$}{Real Euclidean space of dimension $n$}
\nomenclature[A06]{$\mathbb{C}^{n}$}{Complex vector space of dimension $n$}
\nomenclature[A07]{$\mathbb{H}^{n}$}{Upper Euclidean half-space of dimension $n$}
\nomenclature[A08]{$\mathbb{S}^{n}$}{Unit sphere of dimension $n$}
\nomenclature[A09]{$e_{j}$}{The $j$th vector of the standard basis of $\mathbb{R}^{n}$}
\nomenclature[A10]{$B_{d}(x,r)$}{Open ball with center $x$ and radius $r$ for the metric $d$}
\nomenclature[A11]{$\overline{B}_{d}(x,r)$}{Closed ball with center $x$ and radius $r$ for the metric $d$}
\nomenclature[A12]{$A\Subset U$}{The set $A$ is compactly contained in $U$}
\nomenclature[A13]{$\operatorname{supp}(f),\ \operatorname{supp}(\mathbf{f})$}{Support of the function $f$ or the section $\mathbf{f}$}
\nomenclature[A14]{$\mathbf{1}_{A}$}{Indicator function of the set $A$}
\nomenclature[A15]{$\operatorname{span}(A)$}{Vector subspace spanned by $A$}
\nomenclature[A16]{$\displaystyle\sum_{i=1}^{k}a_{i}$}{Sum of the terms $a_i$}
\nomenclature[A17]{$\displaystyle\prod_{i=1}^{k}a_{i}$}{Product of the terms $a_i$}
\nomenclature[A18]{$\displaystyle\coprod_{j\in I}A_{j}$}{Disjoint union of the family $(A_j)_{j\in I}$}
\nomenclature[A19]{$\binom{n}{k}$}{Binomial coefficient}
\nomenclature[A20]{$\alpha$}{Multi-index, in the context of coordinate calculations}
\nomenclature[A21]{$D^{\alpha}u$}{Partial derivative of $u$ associated with the multi-index $\alpha$}

\nomenclature[B01]{$(X,\tau)$}{Topological space $X$ with topology $\tau$}
\nomenclature[B02]{$\overline{A}$}{Closure of the set $A$}
\nomenclature[B03]{$\sim$}{Equivalence relation}
\nomenclature[B04]{$\overline{X}^{\lVert\cdot\rVert}$}{Completion of the normed space $(X,\lVert\cdot\rVert)$}
\nomenclature[B05]{$X'$}{Continuous dual of the topological vector space $X$}
\nomenclature[B06]{$V^{*}$}{Algebraic dual of the vector space $V$}
\nomenclature[B06A]{$\overline V$}{Complex conjugate vector space of $V$}
\nomenclature[B07]{$\operatorname{Hom}(V,W)$}{Space of linear transformations from $V$ to $W$}
\nomenclature[B08]{$\mathcal{L}(X,Y)$}{Space of continuous linear operators from $X$ to $Y$}
\nomenclature[B09]{$\langle\cdot,\cdot\rangle$}{Inner product or duality pairing, according to context}
\nomenclature[B10]{$\lVert\cdot\rVert$}{Norm of the space under consideration}
\nomenclature[B11]{$\lVert\cdot\rVert_{p}$}{Norm of a Lebesgue space}
\nomenclature[B12]{$x_{j}\rightharpoonup x$}{Weak convergence of $x_j$ to $x$}
\nomenclature[B13]{$x_{j}\overset{*}{\rightharpoonup}x$}{Weak-$*$ convergence of $x_j$ to $x$}
\nomenclature[B14]{$\sigma(X,Y)$}{Weak topology on $X$ induced by duality with $Y$}
\nomenclature[B15]{$(X,\mathcal{A},\mu)$}{Measure space with $\sigma$-algebra $\mathcal A$ and measure $\mu$}
\nomenclature[B16]{$\lambda_{n}$}{Lebesgue measure on $\mathbb{R}^{n}$}
\nomenclature[B17]{$L^{p}(X,\mathcal{A},\mu)$}{Lebesgue space on the measure space $(X,\mathcal A,\mu)$}
\nomenclature[B18]{$M_{m\times n}(\mathbb{K})$}{Matrices of size $m\times n$ with entries in $\mathbb K$}
\nomenclature[B19]{$A^{T}$}{Transpose of the matrix $A$}
\nomenclature[B20]{$A^{*}$}{Hermitian adjoint of a matrix; for operators, only when the context specifies the Hilbert inner products}
\nomenclature[B20A]{$T^{\dagger}$}{Bounded Hilbert space adjoint of $T$ when it must be distinguished from the differential formal adjoint}
\nomenclature[B21]{$\mathcal D(A)$}{Domain of the linear operator $A$}
\nomenclature[B22]{$\mathcal G(A)$}{Graph of the linear operator $A$}
\nomenclature[B23]{$M^{\perp}$}{Orthogonal complement of the subset $M$ in a Hilbert space}
\nomenclature[B24]{$\mathcal K(X,Y)$}{Space of compact linear operators from $X$ to $Y$}
\nomenclature[B25]{$\mathcal L^m(X_1,\ldots,X_m;Y)$}{Space of continuous $m$-linear operators from $X_1\times\cdots\times X_m$ to $Y$}
\nomenclature[B26]{$\mathcal L_s^m(X;Y)$}{Subspace of continuous symmetric $m$-linear operators}
\nomenclature[B27]{$D^mf(x)$}{Fréchet derivative of order $m$ of $f$ at $x$}
\nomenclature[B28]{$\delta^mf(x)$}{Gâteaux derivative of order $m$ of $f$ at $x$}
\nomenclature[B29]{$\operatorname{Car}(\Omega\times\mathbb R^m,\mathbb R^\ell)$}{Carathéodory functions with the indicated domains and codomains}
\nomenclature[B30]{$\delta_jf(x)$}{Partial Gâteaux derivative of $f$ with respect to the $j$th variable}
\nomenclature[B31]{$D_jf(x)$}{Partial Fréchet derivative of $f$ with respect to the $j$th variable}
\nomenclature[B32]{$\mathcal N_f$}{Nemytski or superposition operator associated with $f$}
\nomenclature[B33]{$\operatorname{Iso}(X,Y)$}{Set of continuous linear isomorphisms from $X$ onto $Y$}
\nomenclature[B34]{$(X_0,X_1)_{\theta,q}$}{Real interpolation space of the compatible couple $(X_0,X_1)$}
\nomenclature[B35]{$[X_0,X_1]_{\theta}$}{Complex interpolation space of the compatible couple $(X_0,X_1)$}
\nomenclature[B36]{$K(t,u;X_0,X_1)$}{Lions--Peetre $K$-functional}
\nomenclature[B37]{$J(t,u;X_0,X_1)$}{Lions--Peetre $J$-functional}
\nomenclature[B38]{$\mathcal S(\mathbb R^n)$}{Schwartz class on $\mathbb R^n$}
\nomenclature[B39]{$\mathcal S'(\mathbb R^n)$}{Space of tempered distributions}
\nomenclature[B40]{$\mathcal F u=\widehat u$}{Fourier transform of $u$ with unitary normalization}
\nomenclature[B41]{$X/Y$}{Quotient of the normed space $X$ by the closed subspace $Y$}
\nomenclature[B42]{$BC^j(U,F),\ BC^\infty(U,F)$}{Functions of class $C^j$ with bounded derivatives up to order $j$; the infinite-order class is the intersection over all orders}
\nomenclature[B43]{$BUC^j(U,F),\ BUC^\infty(U,F)$}{Subspaces of $BC^j$ whose derivatives up to order $j$ are uniformly continuous}

\nomenclature[C01]{$C^{k}(M,N)$}{Maps of class $C^{k}$ between the smooth manifolds $M$ and $N$}
\nomenclature[C02]{$C^{k}_{c}(M)$}{Functions of class $C^{k}$ with compact support in $M$}
\nomenclature[C03]{$C^{\infty}(M)$}{Smooth real-valued functions on $M$}
\nomenclature[C04]{$TM$}{Tangent bundle of the smooth manifold $M$}
\nomenclature[C05]{$T^{*}M$}{Cotangent bundle of the smooth manifold $M$}
\nomenclature[C06]{$T_{p}M$}{Tangent space of $M$ at the point $p$}
\nomenclature[C07]{$\mathfrak{X}(M)$}{Space of smooth vector fields on $M$}
\nomenclature[C08]{$\partial_i$}{Coordinate vector field $\frac{\partial}{\partial} x^i$}
\nomenclature[C09]{$\mathbf{d}x^{i}$}{Coordinate covector field}
\nomenclature[C10]{$df$}{Differential of the function $f$}
\nomenclature[C11]{$T^{(k,l)}(TM)$}{Bundle of $k$-contravariant and $l$-covariant tensors}
\nomenclature[C12]{$F^{i_{1}\dots i_{k}}_{j_{1}\dots j_{l}}$}{Coordinate components of a tensor of type $(k,l)$}
\nomenclature[C13]{$\otimes$}{Tensor product}
\nomenclature[C14]{$\wedge$}{Wedge product}
\nomenclature[C15]{$\bigwedge^{k}(T^{*}M)$}{Bundle of alternating $k$-tensors}
\nomenclature[C16]{$\Omega^{k}(M)$}{Space of differential $k$-forms on $M$}
\nomenclature[C17]{$d\omega$}{Exterior derivative of the differential form $\omega$}
\nomenclature[C18]{$F\restriction_{A}$}{Restriction of $F$ to the subset $A$}
\nomenclature[C19]{$\partial M$}{Boundary of a manifold with boundary $M$}

\nomenclature[D01]{$(M,\mathbf{g})$}{Riemannian manifold $M$ with metric $\mathbf{g}$}
\nomenclature[D02]{$g_{ij}$}{Local coefficients of the Riemannian metric $\mathbf{g}$}
\nomenclature[D03]{$g^{ij}$}{Coefficients of the inverse matrix of $(g_{ij})$}
\nomenclature[D04]{$\langle\cdot,\cdot\rangle_{\mathbf{g}}$}{Inner product induced by the metric $\mathbf{g}$}
\nomenclature[D05]{$\lvert\cdot\rvert_{\mathbf{g}}$}{Pointwise norm induced by the metric $\mathbf{g}$}
\nomenclature[D06]{$\lambda_{\mathbf{g}}$}{Riemann--Lebesgue measure associated with $\mathbf{g}$}
\nomenclature[D07]{$d\lambda_{\mathbf{g}}$}{Integration element with respect to the Riemann--Lebesgue measure}
\nomenclature[D08]{$dV_{\mathbf{g}}$}{Volume form of $\mathbf{g}$ when an orientation has been fixed}
\nomenclature[D09]{$d_{\mathbf{g}}$}{Riemannian distance}
\nomenclature[D10]{$L_{\mathbf{g}}(\gamma)$}{Riemannian length of the curve $\gamma$}
\nomenclature[D11]{$\operatorname{inj}(M,\mathbf{g})$}{Injectivity radius of the Riemannian manifold}
\nomenclature[D12]{$\operatorname{grad}_{\mathbf{g}}f$}{Riemannian gradient of $f$}
\nomenclature[D13]{$\operatorname{div}_{\mathbf{g}}X$}{Riemannian divergence of the vector field $X$}
\nomenclature[D14]{$\Delta_{\mathbf{g}}$}{Laplace--Beltrami operator, $\operatorname{div}_{\mathbf{g}}\operatorname{grad}_{\mathbf{g}}$}
\nomenclature[D15]{$\operatorname{Rm}_{\mathbf{g}}$}{Riemann curvature tensor of the metric $\mathbf{g}$}
\nomenclature[D16]{$\operatorname{Ric}_{\mathbf{g}}$}{Ricci tensor of the metric $\mathbf{g}$}
\nomenclature[D17]{$R=R_{\mathbf{g}}$}{Scalar curvature of the metric $\mathbf{g}$}

\nomenclature[E01]{$\mathbf{E}\longrightarrow M$}{Vector bundle $\mathbf{E}$ over the manifold $M$}
\nomenclature[E02]{$\mathbf{E}^{*}$}{Finite-rank dual vector bundle of $\mathbf{E}$, in accordance with the finite-dimensional convention}
\nomenclature[E02A]{$\overline{\mathbf{E}}$}{Complex conjugate vector bundle of $\mathbf{E}$}
\nomenclature[E03]{$\Gamma(\mathbf{E})$}{Space of smooth sections of $\mathbf{E}$}
\nomenclature[E04]{$\Gamma^{k}(\mathbf{E})$}{Space of sections of class $C^{k}$ of $\mathbf{E}$}
\nomenclature[E05]{$\Gamma^{k}_{c}(\mathbf{E})$}{Compactly supported sections of class $C^{k}$ of $\mathbf{E}$}
\nomenclature[E06]{$\Gamma_{c}(\mathbf{E})$}{Smooth compactly supported sections of $\mathbf{E}$}
\nomenclature[E07]{$\mathbf{h}_{\mathbf{E}}$}{Bundle metric on the vector bundle $\mathbf{E}$}
\nomenclature[E07A]{$\boldsymbol{\mathcal{R}}_{\mathbf{E}}\colon\overline{\mathbf{E}}\to \mathbf{E}^*$}{Complex-linear Riesz isomorphism, $\boldsymbol{\mathcal{R}}_{\mathbf{E}}(\overline v)(w)=\mathbf{h}_{\mathbf{E}}(w,v)$}
\nomenclature[E08]{$\langle\cdot,\cdot\rangle_{\mathbf{h}_{\mathbf{E}}}$}{Pointwise inner product induced by $\mathbf{h}_{\mathbf{E}}$}
\nomenclature[E09]{$\lvert\cdot\rvert_{\mathbf{h}_{\mathbf{E}}}$}{Pointwise norm induced by $\mathbf{h}_{\mathbf{E}}$}
\nomenclature[E10]{$\nabla^{\mathbf{E}}$}{Connection on the vector bundle $\mathbf{E}$}
\nomenclature[E11]{$\nabla_{i}$}{Covariant derivative in the direction of the coordinate vector field $\boldsymbol{\partial}_i$}
\nomenclature[E12]{$\nabla^{k}\mathbf{F}$}{The $k$th covariant derivative of the tensor field or section $\mathbf{F}$}
\nomenclature[E13]{$\nabla^{k}_{w}\mathbf{F}$}{The $k$th weak covariant derivative of $\mathbf{F}$}
\nomenclature[E14]{$\boldsymbol{\mathcal{D}}(M,\mathbf{E})$}{Space of test sections of the bundle $\mathbf{E}$}
\nomenclature[E15]{$\boldsymbol{\mathcal{D}}'(M,\mathbf{E},W)$}{$W$-valued distributions in the vector bundle $\mathbf{E}$}
\nomenclature[E16]{$\boldsymbol{\mathcal{E}}'$}{Continuous dual bundle of a Banach vector bundle $\boldsymbol{\mathcal{E}}$}
\nomenclature[E17]{$V_e\mathbf{E},\ V\mathbf{E}$}{Vertical space at $e\in \mathbf{E}$ and vertical subbundle of $T\mathbf{E}$}
\nomenclature[E18]{$w_e^{\mathrm V}$}{Vertical lift of $w\in \mathbf{E}_{\pi_{\mathbf{E}}(e)}$ at $e$}
\nomenclature[E19]{$K^{\nabla^{\mathbf{E}}}\colon T\mathbf{E}\longrightarrow \mathbf{E}$}{Connection map induced by $\nabla^{\mathbf{E}}$}
\nomenclature[E20]{$H_e^{\nabla^{\mathbf{E}}}\mathbf{E},\ H^{\nabla^{\mathbf{E}}}\mathbf{E}$}{Horizontal space at $e$ and horizontal subbundle induced by $\nabla^{\mathbf{E}}$}
\nomenclature[E21]{$X_e^{\mathrm H}$}{Horizontal lift of $X\in T_{\pi_{\mathbf{E}}(e)}M$ at $e$}
\nomenclature[E22]{$\pi^*\nabla^{\mathbf{E}}$}{Pullback connection on $\pi^*\mathbf{E}$}
\nomenclature[E23]{$\nabla^{\mathrm H},\ \nabla^{\mathrm V}$}{Horizontal and vertical covariant derivatives on the cotangent bundle}

\nomenclature[F01]{$L^{p}(\Omega)$}{Lebesgue space on the open subset $\Omega\subseteq\mathbb{R}^{n}$}
\nomenclature[F02]{$L^{p}(M)$}{Lebesgue space with respect to the measure $\lambda_{\mathbf{g}}$}
\nomenclature[F03]{$L^{p}(M,\mathbf{E})$}{Space of $p$-integrable sections of the bundle $\mathbf{E}\longrightarrow M$}
\nomenclature[F04]{$L^{p}(M,\mathbb{R}^{n})$}{Space of measurable functions from $M$ to $\mathbb{R}^{n}$ whose Euclidean bundle norm is $p$-integrable}
\nomenclature[F05]{$W^{m,p}(\Omega)$}{Sobolev space of order $m$ on $\Omega$}
\nomenclature[F06]{$W^{m,p}(M)$}{Sobolev space of order $m$ on the manifold $M$}
\nomenclature[F07]{$W^{m,p}(M,\mathbf{E})$}{Sobolev space of sections of the bundle $\mathbf{E}$}
\nomenclature[F07A]{$W^{m,p}_{\mathrm{loc}}(M,\mathbf{E})$}{Local Sobolev space of sections of the bundle $\mathbf{E}$}
\nomenclature[F08]{$W^{m,p}_{0}(M)$}{Closure of $C^{\infty}_{c}(\operatorname{Int}M)$ in $W^{m,p}(M)$}
\nomenclature[F09]{$W^{m,p}_{0}(M,\mathbf{E})$}{Closure of $\Gamma_{c,\operatorname{Int}(M)}(\mathbf{E})$ in $W^{m,p}(M,\mathbf{E})$}
\nomenclature[F10]{$H^{s,p}(M)$}{Bessel potential space of order $s$ on the manifold $M$}
\nomenclature[F10A]{$H^{s,p}(M,\mathbf{E})$}{Bessel Sobolev space of sections, defined by uniform localization}
\nomenclature[F10B]{$H^s_{\mathrm{loc}}(M,\mathbf{E})$}{Local Sobolev space of sections of the bundle $\mathbf{E}$}
\nomenclature[F10C]{$H_{L}^{s,p}(M)$}{Intrinsic Bessel scale associated with $L=\Delta_B=-\Delta$}
\nomenclature[F10D]{$H_{L_{\mathbf{E}}}^{s,p}(M,\mathbf{E})$}{Intrinsic Bessel scale associated with the connection Laplacian $L_{\mathbf{E}}$}
\nomenclature[F11]{$p^{*}$}{Critical Sobolev exponent corresponding to the context}
\nomenclature[F12]{$\operatorname{Tr},\ \boldsymbol{\operatorname{Tr}}^{\mathbf{E}}$}{Trace operators for functions and for sections of the bundle $\mathbf{E}$, respectively}
\nomenclature[F13]{$H^{s,p}(\mathbb R^n)$}{Bessel potential space of order $s$}
\nomenclature[F14]{$W^{s,p}(\Omega)$}{Sobolev--Slobodeckij space of noninteger order $s$}
\nomenclature[F14A]{$W^{s,p}(M,\mathbf{E})$}{Intrinsic Sobolev--Slobodeckij space of sections of the bundle $\mathbf{E}$, defined using covariant derivatives and parallel transport}
\nomenclature[F15]{$[u]_{W^{s,p}(\Omega)}$}{Gagliardo--Slobodeckij seminorm}
\nomenclature[F16]{$B^s_{p,q}(M),\ B^s_{p,q}(M,\mathbf{E})$}{Besov spaces of functions and sections on $M$}
\nomenclature[F17]{$F^s_{p,q}(M),\ F^s_{p,q}(M,\mathbf{E})$}{Triebel--Lizorkin spaces of functions and sections on $M$}
\nomenclature[F18]{$J^s=(I-\Delta)^{\frac{s}{2}}$}{Bessel operator of order $s$}

\nomenclature[G01]{$\mathbf{PDO}^{(m)}(\mathbf{E},\mathbf{F})$}{Partial differential operators of order at most $m$ from $\mathbf{E}$ to $\mathbf{F}$}
\nomenclature[G02]{$P_h^{*}$}{Hermitian formal adjoint of the differential operator $P$, with respect to the bundle metrics and the Riemannian measure}
\nomenclature[G02A]{$P'$}{Complex-bilinear formal transpose, $P'\colon\Gamma(\mathbf{F}^*)\to\Gamma(\mathbf{E}^*)$}
\nomenclature[G03]{$\boldsymbol{\sigma}_{m}(P)$}{Principal symbol of order $m$ of the operator $P$}
\nomenclature[G04]{$d^{\mathbf{E}}$}{Exterior derivative with values in the bundle $\mathbf{E}$}
\nomenclature[G05]{$\delta^{\mathbf{E}}$}{Hermitian formal adjoint of $d^{\mathbf{E}}$}
\nomenclature[G06]{$\Delta_{B}$}{Bochner Laplacian with the positive sign convention, $(\nabla)_h^{*}\nabla$}
\nomenclature[G07]{$\Delta_{H}$}{Hodge Laplacian}
\nomenclature[G08]{$\Delta_{p,g}$}{$p$-Laplacian associated with the metric $g$}
\nomenclature[G09]{$S^{m}_{1,0}$}{Hörmander symbol class of order $m$ and type $(1,0)$}
\nomenclature[G10]{$S^{m}_{\mathrm{cl}}$}{Class of classical symbols of order $m$}
\nomenclature[G11]{$\Psi^{m}(M;\mathbf{E},\mathbf{F})$}{Pseudodifferential operators of order $m$ from sections of $\mathbf{E}$ to sections of $\mathbf{F}$}
\nomenclature[G12]{$\operatorname{WF}(u)$}{Wavefront set of the distribution $u$}
\nomenclature[G13]{$S^{-\infty}(U\times\mathbb R^n;\operatorname{Hom}(\mathbb C^r,\mathbb C^s))$}{Euclidean smoothing symbols}
\nomenclature[G14]{$\Psi^{-\infty}(M;\mathbf{E},\mathbf{F})$}{Properly supported smoothing operators on $M$}
\nomenclature[G15]{$\Psi^m_{\mathrm u,\rho}(M;\mathbf{E},\mathbf{F})$}{Radius-$\rho$ stage of the uniform pseudodifferential calculus}
\nomenclature[G16]{$\Psi^{-\infty}_{\mathrm u}(M;\mathbf{E},\mathbf{F})$}{Uniformly smoothing operators of finite propagation}

\nomenclature[H01]{$X,Y$}{Banach spaces serving as model spaces}
\nomenclature[H02]{$DI(u)$}{Fréchet derivative of the functional $I$ at $u$}
\nomenclature[H03]{$T'M$}{Continuous cotangent bundle of a Banach manifold}
\nomenclature[H04]{$dF_{p}$}{Differential of a map between Banach manifolds}
\nomenclature[H05]{$T_{p}^{\operatorname{kin}}M$}{Kinematic tangent space of the Banach manifold $M$}
\nomenclature[H06]{$T_{p}^{\operatorname{op}}M$}{Operational tangent space of the Banach manifold $M$}
\nomenclature[H07]{$T_{p}N$}{Tangent space of the Banach submanifold $N$ at $p$}
\nomenclature[H08]{$\operatorname{GL}(E)$}{Group of continuous invertible linear operators on the Banach space $E$}
\nomenclature[H09]{$\flat_{\mathbf{g}},\ \sharp_{\mathbf{g}}$}{Musical morphisms induced by a strong Riemannian metric}
\nomenclature[H10]{$F_x,\ F_x^*$}{Finsler norm on $T_xM$ and its dual norm on the continuous dual $T'_xM$}
\nomenclature[H11]{$d_F$}{Length distance induced by a Banach--Finsler structure}
\nomenclature[H12]{$W^{r,p}(M,N)$}{Manifold of Sobolev maps from $M$ to $N$}
\nomenclature[H13]{$\mathcal A^r(\mathbf{P})$}{Affine Sobolev space of connections on a principal bundle $\mathbf{P}$}
\nomenclature[H14]{$\mathcal G^{r+1}(\mathbf{P})$}{Gauge group completed in Sobolev regularity}
\nomenclature[H15]{$\operatorname{Diff}^s(M)$}{Sobolev completion of the diffeomorphism group}
\nomenclature[H16]{$H^\infty=\bigcap_{s\geq s_0}H^s$}{Smooth core of a Hilbert scale}
\nomenclature[H17]{$\varprojlim H^s$}{Projective limit of a Hilbert scale}
\nomenclature[H18]{$\operatorname{Met}^r(M)$}{Open space of Riemannian metrics of Sobolev regularity $r$}
\nomenclature[H19]{$C_b^j(M,\mathbf{E})$}{Sections of class $C^j$ whose covariant derivatives up to order $j$ are bounded}
\nomenclature[H20]{$C_b^{\infty,r}(M,N)$}{Smooth maps whose covariant derivatives of $df$ up to order $r-1$ are bounded}
\nomenclature[H21]{$\operatorname{Met}(I,B_k)(M)$}{Complete metrics with positive injectivity radius and curvature bounded up to order $k$}

\nomenclature[I01]{$I$}{Variational functional or energy functional}
\nomenclature[I02]{$DI(u)=0$}{Critical point condition for the functional $I$}
\nomenclature[I03]{$J^{-1}(c)$}{Level set at level $c$ of the functional $J$}
\nomenclature[I04]{$\mathcal{N}$}{Nehari manifold associated with a functional}
\nomenclature[I05]{$I\restriction_{\mathcal N}$}{Restriction of the functional $I$ to the Nehari manifold}
\nomenclature[I06]{$T_{u}\mathcal{N}$}{Tangent space of the Nehari manifold at $u$}
\nomenclature[I07]{$\lambda$}{Lagrange multiplier, when it occurs in a variational constraint}
\nomenclature[I08]{$K_c$}{Set of critical points at level $c$}
\nomenclature[I09]{$I^{\leq a}$}{Closed sublevel set $\{u\mid I(u)\leq a\}$}
\nomenclature[I09A]{$I^{<a}$}{Open sublevel set $\{u\mid I(u)<a\}$}
\nomenclature[I10]{$(PS)_c$}{Palais--Smale condition at level $c$}
\nomenclature[I11]{$\Phi_t$}{Time-$t$ map of a deformation flow}

\printnomenclature
\tableofcontents

\cleardoublepage
\chapter*{Introduction}
\addcontentsline{toc}{chapter}{Introduction}
\markboth{Introduction}{Introduction}

When we begin studying a differential equation on a manifold, something familiar happens: we choose a chart and, for a moment, the problem becomes one defined on an open subset of $\mathbb R^n$. There we can differentiate, integrate, and use many of the tools of Euclidean analysis. The difficulty arises when we want to leave that chart. We must then know which parts of the calculation depended on coordinates, which quantities have intrinsic meaning, and how to combine the estimates obtained locally to draw a conclusion about the entire manifold.

Global analysis studies precisely this passage from the local to the global. Coordinates do not disappear: they remain an indispensable computational tool. What matters is understanding which information survives a change of chart and how geometry modifies arguments we already know in $\mathbb R^n$. A connection allows us to differentiate sections of a bundle; a metric allows us to measure and integrate them; Sobolev spaces make it possible to work with solutions that are not smooth; spectral theory and heat relate local operators to global information; and variational methods turn many equations into energy problems on infinite-dimensional spaces.

This is the viewpoint that guides the book. Rather than presenting differential geometry, functional analysis, and partial differential equations separately, we want to see how these theories meet in the study of a concrete problem on a manifold. Throughout the text, we will therefore return to the same pattern: first formulate the object intrinsically, then express it in coordinates when a calculation is needed, and finally check which hypotheses allow us to pass from local estimates back to a global statement.

Much of this book grew out of teaching in the graduate program in Mathematical Sciences, in the courses \emph{Nonlinear analysis on Riemannian manifolds I and II} and \emph{Variational methods on Riemannian manifolds I and II}. These courses were led by Dr.~María de los Ángeles Sandoval Romero, who taught the nonlinear analysis courses jointly with Carlos Daniel Velázquez Mendoza, M.Sc., and the variational courses jointly with Dr.~Romulo Diaz Carlos.

In the classroom, many arguments can be completed with an explanation at the blackboard, an example, or a reference to something covered a few weeks earlier. Turning that material into a book required these connections to be put into writing. This is the source of our concern to state the hypotheses of each result carefully, recover the intermediate steps of proofs, and show explicitly where geometry enters an argument that might at first appear purely analytic.

\medskip
\noindent\textbf{Note to the reader.}
This manuscript is a work in progress and may still undergo minor corrections
and editorial refinements. Comments, suggestions, and corrections are very welcome.

\section*{From local calculus to intrinsic geometry}

The starting point is differential calculus. Newton and Leibniz developed, in different languages, tools for studying velocities, tangents, areas, and accumulated quantities. Over time, these ideas acquired increasingly precise formulations: first in one variable, then in several variables, and finally in more general spaces. The derivative ceased to be viewed merely as a computational rule and came to be understood as the linear approximation describing the first-order behavior of a nonlinear object.

Differential geometry then required an answer to a new question: what does differentiation mean when the space on which we work is curved? Gauss showed that the curvature of a surface can be recognized from measurements made within the surface itself, without reference to the ambient space. Riemann extended this idea to arbitrary dimensions and proposed describing geometry through a way of measuring lengths and angles infinitesimally. Later, Poincaré brought another difficulty to the forefront: knowing each small region of a space well is not necessarily enough to understand its global shape. The work of Gauss, Riemann, and Poincaré thus marks the passage from calculus on embedded figures to intrinsic geometry, and then to the study of truly global properties \cite{Gauss1828,Riemann1868,Poincare1895}.

\begin{figura}[htbp]
  \begin{minipage}[t]{0.18\textwidth}
    \centering
    \includegraphics[width=\linewidth,height=2.65cm,keepaspectratio]{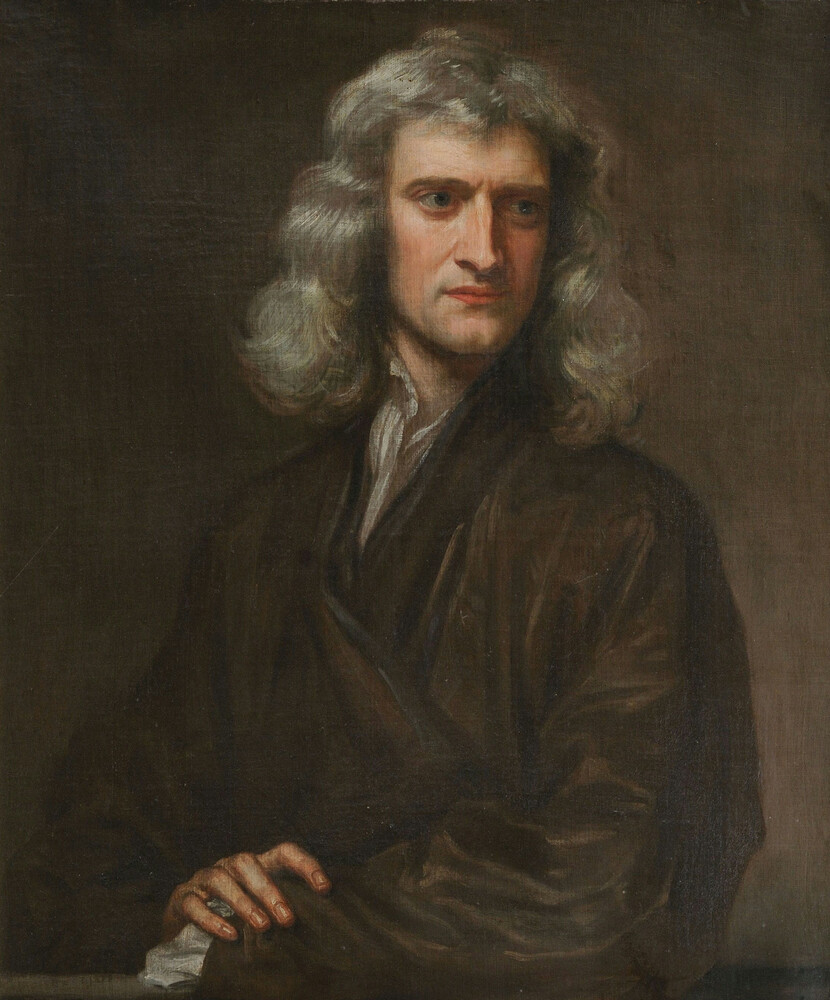}\par
    \smallskip
    \footnotesize Isaac Newton\\(1642--1727)
  \end{minipage}\hfill
  \begin{minipage}[t]{0.18\textwidth}
    \centering
    \includegraphics[width=\linewidth,height=2.65cm,keepaspectratio]{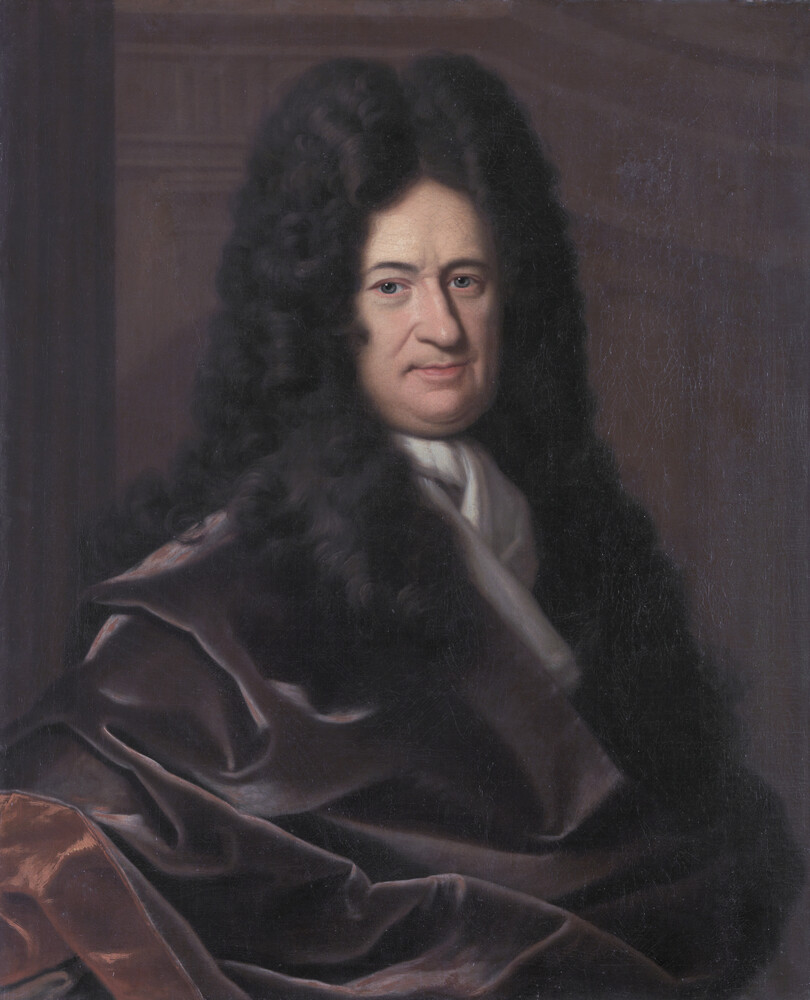}\par
    \smallskip
    \footnotesize G. W. Leibniz\\(1646--1716)
  \end{minipage}\hfill
  \begin{minipage}[t]{0.18\textwidth}
    \centering
    \includegraphics[width=\linewidth,height=2.65cm,keepaspectratio]{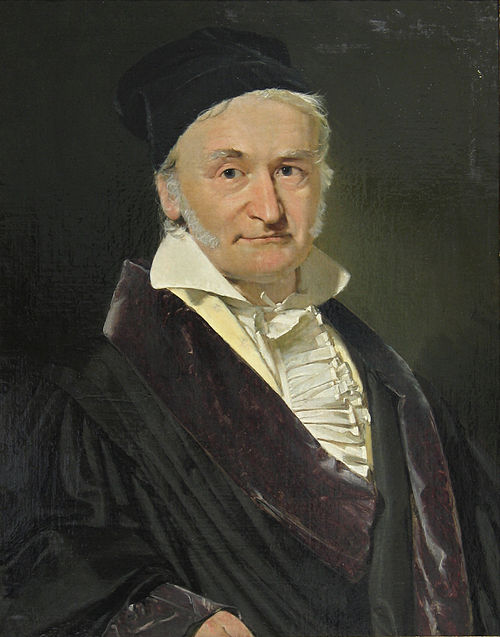}\par
    \smallskip
    \footnotesize C. F. Gauss\\(1777--1855)
  \end{minipage}\hfill
  \begin{minipage}[t]{0.18\textwidth}
    \centering
    \includegraphics[width=\linewidth,height=2.65cm,keepaspectratio]{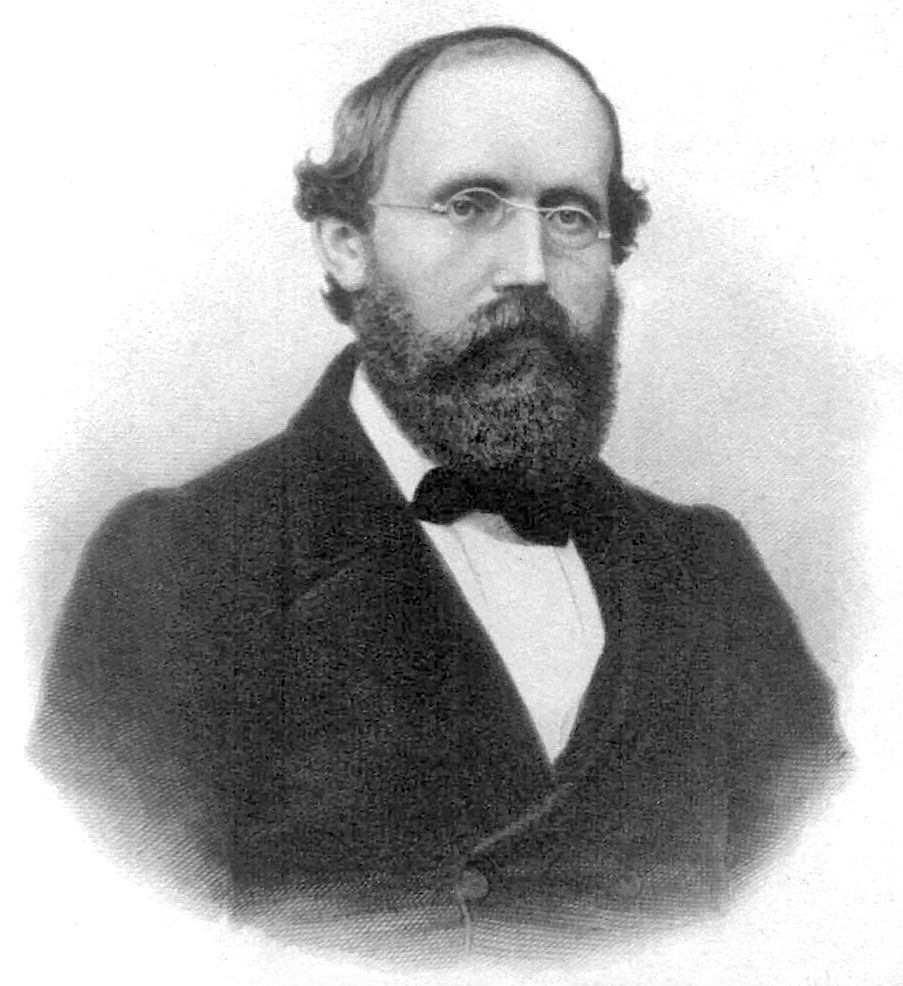}\par
    \smallskip
    \footnotesize B. Riemann\\(1826--1866)
  \end{minipage}\hfill
  \begin{minipage}[t]{0.18\textwidth}
    \centering
    \includegraphics[width=\linewidth,height=2.65cm,keepaspectratio]{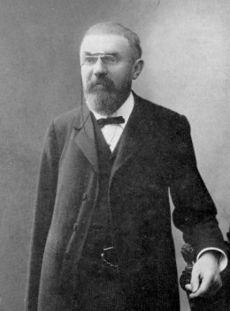}\par
    \smallskip
    \footnotesize H. Poincaré\\(1854--1912)
  \end{minipage}
  \caption[From calculus to global geometry]{From infinitesimal calculus to global geometry. Newton and Leibniz established the fundamental languages of calculus; Gauss discovered the intrinsic nature of surface curvature; Riemann extended metric geometry to arbitrary dimensions; Poincaré made the global organization of spaces a central mathematical problem.}
  \label{fig:introduccion-calculo-geometria-global}
\end{figura}

\begin{semblanzaHistorica}{Gauss, Riemann, and the intrinsic viewpoint}
In his \emph{Disquisitiones generales circa superficies curvas} of 1828, Gauss studied the geometry of a surface through its first fundamental form and established that Gaussian curvature is intrinsic. Riemann's habilitation lecture of 1854, published posthumously in 1868, changed the scale of the problem: a metric could be considered on spaces of any dimension, and curvature became a property of the space itself. The tensor language and connections subsequently developed by Christoffel, Ricci, Levi--Civita, and Cartan gave this vision an operational form. In this book, that history reappears concretely: tangent spaces and tensors are constructed first; the metric, connection, geodesics, and curvature are then introduced; and finally these structures become coefficients, operators, and energies.
\end{semblanzaHistorica}

This historical development also explains the order we will follow. The covariant derivative, Laplacian, Sobolev norms, and principal symbol will be introduced, whenever possible, as intrinsic objects. Coordinates will appear when they are needed for calculations, comparison with a Euclidean model, or the proof of an estimate. We do not seek to avoid local calculations, but to use them without confusing the coordinate expression with the geometric object it represents.

\section*{The emergence of function spaces}

Geometry is not the only difficulty. When the unknown in an equation is a function or a section, the set of possible solutions no longer lies in a finite-dimensional space. To study a minimizing sequence, we need to know in what sense it can converge; to discuss weak derivatives, we need a space that records both size and regularity; to apply variational methods, we need a structure in which continuity, compactness, and differentiation make sense. This was one of the fundamental changes in twentieth-century analysis. In 1906, Fréchet introduced a systematic theory of metric spaces and functional calculus \cite{Frechet1906}, and in 1932 Banach organized a decisive portion of linear analysis around complete normed spaces \cite{Banach1932}. Functions could then be studied as points of a space that, although infinite-dimensional, had an appropriate geometry and topology.

Sobolev showed that a function can have derivatives in a sense compatible with integration by parts even when it is not classically differentiable, and his embedding theorems related the number of integrable derivatives to the regularity of the function \cite{Sobolev1938}. Shortly afterward, Laurent Schwartz developed distribution theory, which allows singular objects such as the Dirac delta to be treated rigorously and extends differentiation to a much broader class \cite{Schwartz1950}. These ideas made it natural to begin with weak solutions, obtain estimates, and then ask whether the equation itself forces regularity to be recovered.

\begin{semblanzaHistorica}{Sobolev and Schwartz: regularity without classical smoothness}
The conceptual change lay not in abandoning derivatives but in recognizing them through their integral effects. In a Sobolev space, the derivative is identified by an integration-by-parts formula; in distribution theory, that same formula becomes the definition of differentiation. This perspective allows us to begin with data of low regularity, obtain a priori estimates, and then prove that the solution has more regularity than initially seemed available. The route we will follow retraces this development: weak derivatives, Sobolev spaces, distributions, Fourier analysis, interpolation, fractional regularity, and, later, their geometric versions.
\end{semblanzaHistorica}

\section*{What makes global analysis global}

In this context, the adjective \emph{global} does not simply mean that the manifold is large or cannot be covered by a single chart. It means that the answer may depend on how the entire space is organized. An elliptic equation is written locally, yet the dimension of its kernel may contain topological information; a Sobolev inequality is proved in charts, yet its constants may depend on volume, curvature, or behavior at infinity; a geodesic satisfies an ordinary differential equation, yet its extendibility for all time is a global property; Ricci flow is a parabolic equation for a metric, yet the study of its singularities ultimately reveals the topology of the manifold.

This relationship between analysis and topology became particularly clear through Hodge theory, Morse theory, and the study of elliptic operators. The Atiyah--Singer theorem subsequently gave that relationship a much more general form: the analytic index of an elliptic operator agrees with an index constructed from its symbol and the topology of the manifold \cite{AtiyahSinger1968I}. Microlocal analysis, the pseudodifferential calculus, and $K$-theory arise naturally on the way to this statement. Instead of presenting them as separate topics, we will develop them as parts of a single chain of ideas.

Another way to obtain global information is to let the geometry evolve. In 1982, Hamilton introduced Ricci flow to deform a metric according to its Ricci curvature \cite{Hamilton1982}. Shortly afterward, DeTurck identified the gauge that reveals a strictly parabolic system behind the diffeomorphism invariance \cite{DeTurck1983}. Two decades later, Perelman introduced monotone quantities, the noncollapsing principle, and new tools for studying singularities and completing the geometrization program in dimension three \cite{Perelman2002,Perelman2003Surgery,Perelman2003Extinction}. This example captures the spirit of geometric analysis: geometry suggests an equation, analysis makes it possible to control its solutions, and the behavior of those solutions ultimately provides information about the original space.

\begin{semblanzaHistorica}{Atiyah--Singer: a bridge between two ways of counting}
For a Fredholm operator, the index is the difference between the dimension of the solution space and the dimension of the obstruction space. At first sight, it is an analytic number: it is obtained by solving an equation. The Atiyah--Singer theorem asserts that this same integer can be calculated from topological data of the principal symbol. The result does not erase the distance between analysis and topology; it explains why they describe the same phenomenon from different ends. For this reason, the part devoted to index theory does not begin with the final formula: it first constructs the pseudodifferential calculus, Fredholm theory, the symbol class, Bott periodicity, and the Thom isomorphism.
\end{semblanzaHistorica}

\section*{A contemporary perspective}

The preceding questions remain open in many directions and today constitute a broad part of global and geometric analysis. Elliptic and spectral theory study how geometry influences the kernel of an operator, its spectrum, and heat. On noncompact manifolds, bounded geometry allows local estimates to be recovered with constants independent of the point. In the presence of a boundary, traces, boundary terms, and the need for precise control of the normal direction arise. Spaces of maps, connections, and metrics lead naturally to infinite-dimensional geometry, while geometric flows use evolution, rescaling, and compactness to study curvature and singularities.

Recent problems directly related to the themes of the book include new forms of approximation by smooth objects in Sobolev spaces on manifolds \cite{chan2024meyersserrin}; the development of Sobolev, Besov, and Bessel potential scales on manifolds with boundary, bounded geometry, or singularities \cite{Amann2025FunctionSpaces}; characterizations of traces that incorporate the connection and curvature \cite{VanSchaftingenWinter2025}; and intrinsic formulations of integration by parts and Sobolev spaces on vector bundles \cite{daniel2026geometricintegrationpartssobolev}. These examples are different, but they share a common concern: determining precisely which geometric hypotheses allow a familiar analytic tool to be transferred to the global setting.

For this reason, throughout the book we will not treat settings requiring different arguments as interchangeable. Compactness can replace uniform control at infinity; a boundary introduces traces and boundary terms; the curvature of a connection reappears when derivatives are commuted; and a weak metric in infinite dimensions does not automatically possess the properties of a strong metric. Distinguishing these hypotheses is more than a matter of precision in statements: it is often exactly what makes it possible to understand why a proof works.

\section*{The approach to exposition}

The book is written so that, whenever possible, the reader can follow the full development of a tool from the Euclidean case to its geometric version. First, the local result is formulated with the hypotheses it actually needs; we then study what happens under multiplication by cutoff functions, changes of coordinates, or passage from one trivialization to another; finally, the result is globalized. Vector bundles also bring bundle metrics and connections into play; noncompact geometry requires uniform control of constants; and in the presence of a boundary, adapted charts and Fermi coordinates allow the normal direction to be separated from the tangential directions.

For the same reason, many proofs are developed in greater detail than would be necessary in a text intended solely as a reference. When an argument is localized, the cutoffs, supports, and norms being compared are specified; when a partition of unity is used, the contribution of its terms is shown; and when a proof depends on an external result, the point at which it is used is identified. The aim is for the logical chain to be reconstructible from the text itself, so that phrases such as ``the argument is local'' or ``this follows by a partition of unity'' do not conceal the steps where the real substance lies.

The longer chapters can be read at different levels of depth. A first reading allows the reader to follow the motivation, definitions, and main results; a closer reading allows the proofs to be reconstructed and the roles of the metric, connection, compactness, and boundary to be identified. The examples and exercises then help verify that the constructions are more than formal. Rather than memorizing a list of theorems, the goal is to recognize a few ideas that recur throughout the book in different forms.

\section*{Organization of the book}

The first part assembles the geometric language used throughout the book. It begins with smooth manifolds, maps, tangent spaces, submanifolds, vector bundles, flows, differential forms, and Lie groups and algebras. Riemannian metrics, principal and associated bundles, connections, parallel transport, geodesics, distance, curvature, and integration then follow. These tools also allow the construction of Clifford algebras, ${\rm Spin}$ and ${\rm Spin}^{c}$ structures, spinor bundles, and their connections. Dirac operators are introduced later, once a general theory of differential operators on bundles is available. The general references for this part are \cite{LeeS,LeeR,Nicolaescu2006,BleeckerBoossIndex}.

The second part serves as a Euclidean laboratory. Integer-order Sobolev spaces are constructed from weak derivatives, distributions and Fourier analysis are developed, and interpolation is introduced. Bessel potentials and Slobodeckij, Besov, and Triebel--Lizorkin spaces then organize regularity of real order. This is the setting in which we prove the multiplication, localization, change-of-variables, embedding, compactness, and trace results needed for the passage to manifolds \cite{Adams,Evans,Leoni2017,Leoni2023,BerghLofstrom1976,Triebel1,Triebel2,Schneider2021}.

With these tools available, the third part turns to analysis on manifolds and vector bundles. Sobolev spaces are defined through weak covariant derivatives, and differential operators and bundle-valued distributions are developed. The integration-by-parts formula, Kato and Poincaré inequalities, embeddings, compactness, and trace theorems appear here. The connection Laplacian and the heat semigroup then provide an intrinsic construction of spaces of fractional order. Finally, bounded geometry allows these tools to be retained on noncompact manifolds with uniform estimates, and traces on submanifolds and boundaries to be studied \cite{Aubin,Hebey1,Hebey2,Eichhorn2007,Nicolaescu2006,grigor2009heat,traces,Amann2025FunctionSpaces}.

The fourth part is devoted to index theory. The pseudodifferential calculus allows parametrices for elliptic operators to be constructed and their regularity to be studied; Fredholm theory turns this information into an analytic index; $K$-theory organizes symbol classes; and Bott periodicity together with the Thom isomorphism leads to the topological index. These constructions finally come together in the Atiyah--Singer theorem for closed manifolds and several of its fundamental geometric examples.

The fifth part concerns geometric evolution equations. We begin with a parabolic theory for quasilinear systems on closed manifolds and apply it to Ricci flow. Diffeomorphism invariance leads to DeTurck's procedure; curvature evolution equations, maximum principles, a priori estimates, and rescalings prepare the study of singularities. The part concludes with Hamilton compactness and Perelman's tools: monotone quantities, noncollapsing, $\kappa$-solutions, and surgery.

The sixth part changes scale and takes geometry into infinite dimensions. Fréchet calculus on Banach spaces leads to Banach manifolds, bundles, connections, flows, and Hilbert or Banach--Finsler structures. Spaces of maps, metrics, and connections, gauge groups, and diffeomorphism groups then enter the picture, with fixed regularity carefully distinguished from ILH scales. The application chapters show how harmonic maps, Yang--Mills theory, and the Euler equations can be viewed as geometric problems on configuration spaces \cite{brezis_functional_2011,Carothers_2004,Lang1999,Eells1966,EbinMarsden1970}.

The seventh part brings together the variational methods we will use to prove the existence of solutions. Weak formulations and Nemytski operators allow energy functionals to be constructed; the direct method combines coercivity, reflexivity, and semicontinuity; the Palais--Smale condition and pseudogradients allow sublevel sets to be deformed; and regular constraints and the Nehari method allow the search for nontrivial critical points. The emphasis is on how these techniques adapt to the function space and the geometry of the problem, rather than on treating them as isolated procedures \cite{Struwe2008,Palais1970,DrabekMilota2013,AttouchButtazzoMichaille2014,kesavan_functional_2023}.

The appendices collect the tools of topology, mathematical analysis, linear algebra, and tensor notation consulted throughout the book. Their purpose is to let the main development remain continuous without concealing the auxiliary results that support the proofs.

\section*{Reading paths}

Although the chapters are arranged for continuous reading, it is not necessary to work through the entire book before reaching a particular topic. For Sobolev spaces and variational problems on compact manifolds, a natural route passes through Chapters~1--9 and continues with Chapters~29--32. For trace theory and bounded geometry, it is useful to add Chapters~12--14. Index theory uses the Clifford geometry of the first part, the differential construction of Dirac operators in Chapter~8, and Chapters~11 and~15--19. Ricci flow requires Riemannian geometry and elliptic and Sobolev theory, and is developed in Chapters~20--22. Infinite-dimensional geometry can be studied starting in Chapter~23, although its applications frequently use the Sobolev spaces constructed earlier.

These paths continually intersect. Integration by parts reappears in the definition of weak derivatives, formal adjoints, and Euler--Lagrange equations; compactness underlies Rellich--Kondrashov, Fredholm theory, and Palais--Smale; a connection serves both to differentiate sections and to construct Laplacians, study heat, or describe configuration spaces; and the principal symbol begins as a local object and ends by relating regularity, parametrices, and the topological index. An important part of the book consists precisely in making these connections visible.

The book is intended for advanced undergraduate and graduate students, as well as readers approaching global analysis from geometry or differential equations. Basic knowledge of real analysis, topology, linear algebra, and multivariable calculus is assumed; when an additional tool is needed, it will be recalled in the text or a reference will be given to the appropriate appendix. The aim is that, as readers progress through these pages, they will be able to move naturally between three levels that recur throughout: the equation to be studied, the function space in which it is formulated, and the geometry of the space on which it is defined.

\cleardoublepage
\pagenumbering{arabic}

\part{Riemannian manifolds}
\chapter*{Introduction to Part I}
\addcontentsline{toc}{chapter}{Introduction to Part I}
\markboth{Part I. Riemannian manifolds}{Introduction to Part I}

Before studying equations, Sobolev spaces, or differential operators on a manifold, we must address a more elementary question: how can we do calculus on a space that generally has no single coordinate system? Locally, a manifold resembles an open subset of $\mathbb R^n$, so many constructions can be expressed in coordinates. Yet the geometric object is not the chart we choose to describe it. The first aim of this part is to learn how to pass from a local description to a formulation independent of that choice.

Chapter~1 constructs the smooth structure and the objects that arise from it: tangent and cotangent spaces, submanifolds, vector bundles, flows, differential forms, Lie groups, and smooth actions. Chapter~2 adds the Riemannian structure. Connections, parallel transport, geodesics, distance, curvature, Riemannian measure, and the basic differential operators then enter the picture. The part concludes by introducing Clifford algebras, the groups ${\rm Spin}$ and ${\rm Spin}^{c}$, and spinor bundles. Dirac operators will be studied later, once the general theory of differential operators on bundles is available.

These constructions are linked by very concrete needs. Charts allow us to calculate, but the tangent space identifies the directions that exist independently of coordinates. A vector bundle brings together data that vary from point to point; a connection makes it possible to compare and differentiate them; the metric measures lengths, angles, and energies; curvature records the discrepancy that arises when transporting and differentiating in different directions. The intention in this part is for these ideas to emerge as natural responses to geometric problems, rather than as an isolated sequence of definitions.

\begin{semblanzaHistorica}{A shared geometric tradition}
Gauss showed that a surface has geometric properties that can be recognized without reference to the space in which it is drawn. Riemann extended this idea to arbitrary dimensions and made the metric the starting point of intrinsic geometry. In the following decades, Christoffel, Ricci, and Levi--Civita developed tensor calculus and the notion of a connection, while Cartan reorganized much of the theory using differential forms and moving frames. The modern language retains all these perspectives. We will therefore use coordinates, tensors, connections, and differential forms as the problem requires, without treating them as competing languages.
\end{semblanzaHistorica}

  \chapter{Fundamental notions of smooth manifolds}
 \vspace{1em}
\begin{quote}
\emph{“One geometry cannot be more true than another; it can only be more convenient.”} \\
\hfill --- Henri Poincaré, \textit{Science and Hypothesis} (1902)
\end{quote}
\vspace{1em}
Differential geometry begins with a simple observation: a space may look almost Euclidean when viewed up close, yet have a global shape that no single coordinate system can describe. A sphere, a torus, or a surface with boundary admits local calculus, but its most important features emerge only when different regions are compared. The notion of a smooth manifold gives precise form to this interplay between the local and the global.

This first chapter develops the language that will accompany the reader throughout the book. We begin with derivatives in Euclidean spaces, move on to charts and changes of coordinates, and arrive at tangent spaces, vector bundles, submanifolds, flows, and differential forms. Each definition addresses a concrete geometric need and prepares an analytic tool for later use. The exposition follows primarily Lee~\cite{LeeS,LeeR} and draws on Nicolaescu's approach~\cite{Nicolaescu2006}, particularly where smooth geometry begins to provide the foundations for global analysis.
 \section{From calculus to differential geometry}

Differential geometry extends a central idea of calculus: approximating a nonlinear object locally by linear data. Newton and Leibniz independently developed the foundations of differential and integral calculus; Leibniz's notation and Newton's dynamical interpretation continue to be reflected in the modern language of derivatives, curves, and fields.

\begin{figura}[htbp]
    \centering

    \begin{minipage}[t]{0.42\textwidth}
        \centering
        \includegraphics[height=3.4cm]{newton.jpg}

        \smallskip
        \footnotesize Isaac Newton\\
        (1642--1727)
    \end{minipage}
    \hfill
    \begin{minipage}[t]{0.42\textwidth}
        \centering
        \includegraphics[height=3.4cm]{leibniz.jpg}

        \smallskip
        \footnotesize Gottfried Wilhelm Leibniz\\
        (1646--1716)
    \end{minipage}

    \caption[Newton and Leibniz]{Newton and Leibniz developed the foundations
    of differential and integral calculus in different mathematical languages.}
    \label{fig:newton-leibniz}
\end{figura}

The subsequent formalization of calculus led to the modern notions of limit, continuity, and differentiability. On a manifold, coordinates allow us to transfer a problem temporarily to an open subset of $\mathbb R^n$; compatibility between charts ensures that the resulting objects do not depend on this local choice. The derivative is then understood as more than a slope: it is the linear map that approximates a nonlinear map to first order.

Before introducing smooth manifolds, we recall this formulation of the derivative in Euclidean spaces. Fréchet's definition expresses linear approximation precisely and anticipates the version we will later use in Banach spaces.

\begin{definition}\label{def:nociones-fundamentales-norma-conjunto-abierto-funcion-frechet}\index{Frechet differentiability in Euclidean spaces@Fréchet differentiability in Euclidean spaces}\index{Frechet derivative@Fréchet derivative!in Euclidean spaces}
Let $U\subseteq \mathbb{R}^{n}$ be an open set and let
$f\colon U\longrightarrow \mathbb{R}^{m}$ be a function. We say that $f$ is
Fréchet differentiable, or simply differentiable, at $x_{0}\in U$ if there exists
a linear map $L\colon \mathbb{R}^{n}\longrightarrow \mathbb{R}^{m}$
such that $\displaystyle \lim_{h\to 0}
\displaystyle\frac{f(x_0+h)-f(x_0)-L(h)}{\|h\|}=0$. We denote this linear map by $df_{x_{0}}$ or by
$Df(x_{0})$. The map $df_{x_0}$ is unique and is called the derivative
of $f$ at $x_0$. If $f$ is differentiable at every point of $U$, we say that $f$
is differentiable on $U$, or simply differentiable when the domain is
understood.
\end{definition}

In the preceding definition, it is essential that the domain of $f$ be open.
Since $U$ is open, there exists $r>0$ such that $B_{\mathrm{euc}}(x_0,r)\subseteq U$. Therefore, if $h\in\mathbb R^n$ and $\|h\|<r$, then $x_0+h\in U$, and
the linear approximation to $f$ can be tested in every direction in
$\mathbb{R}^n$. This condition allows us to recover a unique
linear map associated with $f$ at the point $x_0$.

The situation changes when the domain is not open. In
that case, the expression $f(x_0+h)$ makes sense only for those $h$ for
which $x_0+h$ belongs to the domain of $f$. Thus the linear approximation
detects the behavior of $f$ only along directions
contained in the domain. As a result, two distinct linear
maps may give the same approximation on the set where
the function is defined.

\begin{example}\label{ej:nociones-fundamentales-punto-transformaciones-lineales-tanto-mismo-modo}
Let $A=\{(x,0)\in \mathbb{R}^2\mid x\in \mathbb{R}\}$ and let
$f\colon A\longrightarrow \mathbb{R}$ be given by $f(x,0)=x$. Consider the point
$0=(0,0)\in A$. Let $L,\widetilde{L}\colon \mathbb{R}^{2}\longrightarrow
\mathbb{R}$ be the linear maps given by $L(x,y)=x$ and
$\widetilde{L}(x,y)=x+y$.

If $h=(h_1,0)\in A$, then
$\displaystyle f(0+h)-f(0)-L(h)=h_1-h_1=0$. Therefore,

\[
\lim_{\substack{h\to 0\\ h\in A}}
\frac{
f(0+h)-f(0)-L(h)
}{
\|h\|
}
=0.
\]

Similarly, if $h=(h_1,0)\in A$, we have
$\displaystyle f(0+h)-f(0)-\widetilde{L}(h)
=h_1-\widetilde{L}(h_1,0)=h_1-h_1=0$. Hence,

\[
\lim_{\substack{h\to 0\\ h\in A}}
\frac{
f(0+h)-f(0)-\widetilde{L}(h)
}{
\|h\|
}
=0.
\]

Thus the linear maps $L$ and $\widetilde{L}$ satisfy the same
approximation condition on $A$, although they are distinct as
linear maps on $\mathbb{R}^2$. The reason is that the set
$A$ allows us to approach the origin only through vectors of the form
$(h_1,0)$. Consequently, the function does not detect the component transverse to the $x$-axis: $L$ and $\widetilde{L}$ agree on $A$ but differ in the transverse directions.
\end{example}

The preceding example explains why the usual definition of the derivative is formulated
for open domains. In differential geometry, however, we need to
consider functions defined on subsets that are not open in
$\mathbb{R}^n$, as happens in the study of manifolds with boundary. In these
cases, a natural way to define regularity is to require
the existence of smooth local extensions.

\begin{definition}[Differentiability by local extension]\label{def:nociones-fundamentales-diferenciabilidad-mediante-extension-local}\index{differentiability by local extension@differentiability by local extension}
Let $A\subseteq \mathbb{R}^n$, let $x_{0}\in A$, and let
$f\colon A\longrightarrow \mathbb{R}^m$. We say that $f$ is differentiable at
$x_0$ by local extension if there exist an open set
$U_{x_0}\subseteq \mathbb{R}^n$, with $x_0\in U_{x_0}$, and a function
$\widetilde{f}\colon U_{x_0}\longrightarrow \mathbb{R}^m$ differentiable at
$x_0$, such that
$\displaystyle \widetilde{f}\restriction_{U_{x_0}\cap A}
=
f\restriction_{U_{x_0}\cap A}$.

If the extension $\widetilde{f}$ can be chosen to be of class $C^k$ in a neighborhood
of $x_0$, we say that $f$ is of class $C^k$ at $x_0$ by local
extension.
\end{definition}

\begin{remark}\label{obs:nociones-fundamentales-definicion-anterior-establece-nocion-regularidad-funciones}
The preceding definition provides a notion of regularity for functions
defined on subsets of $\mathbb{R}^n$, but it does not, by itself,
define an intrinsic derivative of $f$ at $x_0$ as a linear map
$\mathbb{R}^n\to\mathbb{R}^m$. The reason is that the derivative may depend on the
chosen local extension.

The preceding example illustrates this dependence. The function $f\colon A\longrightarrow\mathbb{R}$ admits at least the extensions $\widetilde{f}_1(x,y)=x$ and
$\widetilde{f}_2(x,y)=x+y$. Both agree with $f$ on $A$, but
$\displaystyle d(\widetilde{f}_1)_{(0,0)}(v_1,v_2)=v_1$, whereas
$\displaystyle d(\widetilde{f}_2)_{(0,0)}(v_1,v_2)=v_1+v_2$. Thus the
derivatives of the extensions do not agree. This shows that, for
arbitrary subsets, the existence of smooth local extensions allows us
to define regularity, but does not necessarily yield a well-defined
derivative.
\end{remark}
The derivative is interpreted as the \textit{best linear approximation to the function $f$ at $x_{0}$}: if a function $f\colon U\subseteq \mathbb{R}^{n}\longrightarrow \mathbb{R}$ is differentiable at $x_{0}\in U$, then the graph of the function $H\colon \mathbb{R}^{n}\longrightarrow \mathbb{R}$ given by $H(x)=df_{x_{0}}(x-x_{0})+f(x_{0})$ is the tangent hyperplane to the graph of $f$ at the point $(x_{0},f(x_{0}))$.
\begin{example}\label{ejemplo plano tangente}
 Let $f\colon \mathbb{R}^{2}\longrightarrow \mathbb{R}$ be given by $f(x,y)=x^{2}+xy+y^{2}$. We know that $f$ is differentiable and that $df_{(x,y)}(u,v)=(2x+y)u+(x+2y)v$. The tangent hyperplane to the graph of $f$ at $(x,y,f(x,y))$ is given by $H(u,v)=df_{(x,y)}((u,v)-(x,y))+f(x,y)=(2x+y)(u-x)+(x+2y)(v-y)+x^{2}+xy+y^{2}=2xu + yu + xv + 2yv - x^2 - y^2 - xy$, or, in Cartesian form, is the plane $w=2xu + yu + xv + 2yv - x^2 - y^2 - xy.$
\begin{center}

 \includegraphics[scale=1]{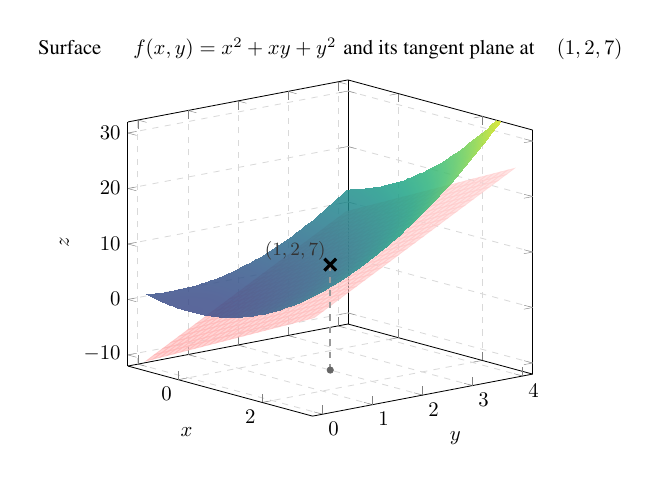}
 \end{center}
\end{example}
The notion of a tangent hyperplane illustrates especially well the power of the \textit{derivative} to provide linear approximations to nonlinear objects, the idea underlying what is known as \textit{nonlinear analysis}.

We can define the tangent hyperplane to a surface if we know a function describing the surface, either as the graph of the function or through a parametrization. However, the notion of a \textit{surface} should not depend on a function describing it: the surface should be a well-defined mathematical object with fixed topological properties. As we will see later, the concept of a \textit{differentiable manifold} resolves these difficulties.

Other approximations to a differentiable function are obtained using polynomial functions. These approximations need not be linear, but they approximate the function increasingly well, thanks to what is known as \textit{Taylor's theorem}. The coefficients of the Taylor expansion are given by higher-order partial derivatives. These derivatives will also help us study the regularity of functions that need not be differentiable, using Sobolev spaces and the concept of a \textit{weak derivative}. Their properties will be essential there, since several properties of weak derivatives are inherited from ordinary partial derivatives. We therefore review these derivatives and the notion of a \textit{\textbf{multi-index}} for readers unfamiliar with this notation.
\subsection{Higher-order partial derivatives}
Before iterating partial derivatives, it is useful to recall the directional derivative. This notion connects Euclidean calculus with a construction that will later appear on manifolds: differentiating a field or a section in the direction of a vector by means of a connection.
\begin{definition}[directional derivative]\label{def:nociones-fundamentales-derivada-direccional}\index{directional derivative}
 Let $U\subseteq \mathbb{R}^{n}$ be open, let $f\colon U\longrightarrow \mathbb{R}^{m}$ be a function, and let $x_{0}\in U$ and $v\in \mathbb{R}^{n}$.
 We say that $f$ is differentiable at $x_{0}$ in the direction of $v$ if
 \[
 \lim_{t\to 0}\frac{f(x_{0}+tv)-f(x_{0})}{t}
 \]
 exists. This limit is called the \textit{\textbf{directional derivative of $f$ in the direction of $v$}} and is denoted by $D_{v}f(x_{0})$.
 The directional derivative for $v=0$ always exists and equals $0$.

 If $f$ is differentiable at every point $x\in U$ in the direction of $v$, we simply say that $f$ has a directional derivative in the direction of $v$ or with respect to $v$.
\end{definition}

These directional derivatives can be iterated to define \textit{higher-order directional derivatives}.
\begin{definition}\label{def:nociones-fundamentales-derivada}\index{higher-order directional derivatives}
Let $f\colon U\subseteq \mathbb{R}^{n}\longrightarrow \mathbb{R}^{m}$ and $x_{0}\in U$.
We define higher-order directional derivatives recursively.

 \begin{enumerate}
 \item If $v_{1},v_{2}\in \mathbb{R}^{n}$ and there exists $r_{1}>0$ such that $B_{\mathrm{euc}}(x_{0},r_{1})\subseteq U$ and $D_{v_{1}}f(x)$ exists for every $x\in B_{\mathrm{euc}}(x_{0},r_{1})$, we say that $f$ has a second directional derivative with respect to $v_{1},v_{2}$ at $x_{0}$ if
 \[
 \lim_{t\to 0}\frac{D_{v_{1}}f(x_{0}+t v_{2}) - D_{v_{1}}f(x_{0})}{t}
 \]
 exists, and we denote this limit by $D_{v_{2},v_{1}}f(x_{0})$.

 \item If $v_{1},\dots,v_{k+1}\in \mathbb{R}^{n}$ and there exists $r_{k}>0$ such that $B_{\mathrm{euc}}(x_{0},r_{k})\subseteq U$ and $D_{v_{k},\dots,v_{1}}f(x)$ is defined and exists for every $x\in B_{\mathrm{euc}}(x_{0},r_{k})$, we say that $f$ has a $(k+1)$st directional derivative with respect to $v_{1},\dots,v_{k+1}$ at $x_{0}$ if
 \[
 \lim_{t\to 0}\frac{D_{v_{k},\dots,v_{1}}f(x_{0}+t v_{k+1}) - D_{v_{k},\dots,v_{1}}f(x_{0})}{t}
 \]
 exists, and we denote this limit by $D_{v_{k+1},\dots,v_{1}}f(x_{0})$.
 \end{enumerate}
 If, for $v_{1},\dots,v_{k}$, $D_{v_{k},\dots,v_{1}}f(x)$ exists for every $x\in U$, we say that $f$ has a directional derivative with respect to $v_{1},\dots,v_{k}$.
\end{definition}

Although the notion of a directional derivative is pointwise, if a function has a directional derivative at every point of its domain in the same direction (or directions), this defines a function with the same domain and codomain as the original function.
\begin{definition}\label{def:nociones-fundamentales-denotamos-funcion-asigna-valor}\index{iterated directional derivatives as functions}
 Let $f\colon U\subseteq \mathbb{R}^{n}\longrightarrow \mathbb{R}^{m}$ and $v_{1},\dots,v_{k}\in \mathbb{R}^{n}$.
 If $D_{v_{k},\dots,v_{1}}f(x)$ exists for every $x\in U$, we denote by
 $D_{v_{k},\dots,v_{1}}f$ the function assigning to each $x\in U$ the value $D_{v_{k},\dots,v_{1}}f(x)$.
\end{definition}

This allows us to define higher-order partial derivatives, which have various uses, including polynomial approximations to sufficiently regular functions.
\begin{definition}\label{def:nociones-fundamentales-orden-de-la-derivada-parcial}\index{order of a partial derivative}
 Let $f\colon U\subseteq \mathbb{R}^{n}\longrightarrow \mathbb{R}^{m}$ and $x_{0}\in U$. Let $\{e_{1},\dots,e_{n}\}$ be the standard basis of $\mathbb{R}^{n}$, and let $i_{1},\dots,i_{k}\in \{1,\dots,n\}$.
 If $\displaystyle D_{e_{i_{k}}\dots e_{i_{1}}}f(x_{0})$ exists, we define the $k$th partial derivative of $f$ with respect to $x^{i_{1}},\dots,x^{i_{k}}$ at $x_{0}$, denoted by $\displaystyle\frac{\partial^{k} f}{\partial x_{i_{k}}\cdots \partial x_{i_{1}}}(x_{0})$, by $\displaystyle\frac{\partial^{k} f}{\partial x_{i_{k}}\cdots \partial x_{i_{1}}}(x_{0})
 :=D_{e_{i_{k}}\dots e_{i_{1}}}f(x_{0})$.

 If $i_{1}=\cdots =i_{k}=i$, we simply denote it by $\displaystyle\frac{\partial^k f}{\partial x_{i}^{k}}(x_{0})$ and call it the $k$th partial derivative of $f$ with respect to $x^{i}$ at $x_{0}$.
 The number $k$ is called the \textit{\textbf{order of the partial derivative}}.

 If $k=0$, we simply write $\displaystyle\frac{\partial^{k} f}{\partial x_{i}^{k}}(x_{0})=f(x_{0})$.

\end{definition}

Multi-index notation provides a compact way to express the existence
of all partial derivatives up to a given order.

\begin{definition}\label{def:nociones-fundamentales-multiindice}\index{multi-index@multi-index}
 An element
 $\alpha=(\alpha_{1},\dots,\alpha_{n})\in\mathbb{N}_{0}^{n}$ is called
 a \textbf{multi-index} of length $n$. Its order is
 $|\alpha|:=\alpha_{1}+\cdots+\alpha_{n}$. If
 $f\colon U\subseteq\mathbb{R}^{n}\longrightarrow\mathbb{R}^{m}$ and the indicated
 derivatives exist, we write
 \[
 D^{\alpha}f
 :=
 \frac{\partial^{|\alpha|}f}
 {\partial x_{1}^{\alpha_{1}}\cdots\partial x_{n}^{\alpha_{n}}},
 \]
 with differentiation understood componentwise when $m>1$, and set
 $D^{0}f:=f$. We denote by $e_j$ the multi-index whose entries are all zero
 except the $j$th, which equals one. We also set
 \[
 \alpha!:=\alpha_1!\cdots\alpha_n!,
 \qquad
 x^\alpha:=(x^1)^{\alpha_1}\cdots(x^n)^{\alpha_n},
 \]
 and write $\beta\leq\alpha$ when
 $\beta_j\leq\alpha_j$ for every $j$.
\end{definition}

The class of functions having continuous partial derivatives up to order
$k$ will be important in defining differentiable structures.

\begin{definition}[functions of class $C^{k}$]\label{def:nociones-fundamentales-funciones-de-clase}\index{function!of class \(C^k\)}
 Let $U\subseteq\mathbb{R}^{n}$ be an open set and let
 $f\colon U\longrightarrow\mathbb{R}^{m}$.
 We say that $f\in C^{0}(U,\mathbb{R}^{m})$ if $f$ is continuous. For
 $k\in\mathbb{N}$, we say that $f\in C^{k}(U,\mathbb{R}^{m})$ if, for
 every $\alpha\in\mathbb{N}_{0}^{n}$ with $|\alpha|\leq k$, the partial
 derivative $D^{\alpha}f$ exists on $U$ and is continuous on $U$. Finally,
 we say that $f\in C^{\infty}(U,\mathbb{R}^{m})$ if
 $f\in C^{k}(U,\mathbb{R}^{m})$ for every $k\in\mathbb{N}_{0}$.

 If $x_{0}\in U$ and $k\in\mathbb{N}_{0}\cup\{\infty\}$, we say that $f$
 is of class $C^{k}$ at $x_{0}$ if there is an open neighborhood
 $V\subseteq U$ of $x_{0}$ such that
 $f\restriction_{V}\in C^{k}(V,\mathbb{R}^{m})$. Functions of class
 $C^{\infty}$ are also called \emph{smooth}.

 We write
 \[
 C^{k}(U,\mathbb{R}^{m})
 :=\{f\colon U\longrightarrow\mathbb{R}^{m}\mid f\text{ is of class }C^{k}\}.
 \]
 If $m=1$, we write $C^{k}(U):=C^{k}(U,\mathbb{R})$.
\end{definition}

The following result, known as the \textit{Clairaut--Schwarz
theorem}, gives local conditions under which
mixed second partial derivatives commute.
\begin{theorem}[Clairaut--Schwarz]\label{clairaut schwarz parciales conmutan}\index{Clairaut Schwarz theorem@Clairaut--Schwarz theorem}
 Let $f\colon U\subseteq \mathbb{R}^{n}\longrightarrow \mathbb{R}^{m}$ and $x_{0}\in U$.
 Suppose that there exists $r>0$ such that
 $B_{\mathrm{euc}}(x_0,r)\subseteq U$. If $i,j\in\{1,\dots,n\}$ and the
 second partial derivatives
 \(\displaystyle\frac{\partial ^{2}f}{\partial x_{i}\partial x_{j}}(x)\)
 and \(\displaystyle\frac{\partial ^{2}f}{\partial x_{j}\partial x_{i}}(x)\)
 exist for every $x\in B_{\mathrm{euc}}(x_{0},r)$ and are continuous at $x_{0}$, then
 \[
 \frac{\partial^{2}f}{\partial x_{i}\partial x_{j}}(x_{0})
 =
 \frac{\partial ^{2}f}{\partial x_{j}\partial x_{i}}(x_{0}).
 \]
\end{theorem}
\begin{proof}
It suffices to work componentwise, so we assume that \(m=1\).
Let \(h,k\ne0\) be small enough that the rectangle with
vertices
\[
x_0,\qquad x_0+h e_i,\qquad x_0+k e_j,\qquad x_0+h e_i+k e_j
\]
is contained in \(B_{\mathrm{euc}}(x_0,r)\), and define
\[
R(h,k):=f(x_0+h e_i+k e_j)-f(x_0+h e_i)
-f(x_0+k e_j)+f(x_0).
\]
We apply the mean value theorem, first in the \(e_i\) direction to
\[
t\longmapsto f(x_0+t e_i+k e_j)-f(x_0+t e_i),
\]
and then in the \(e_j\) direction. There exist
\(\theta_{h,k}\) between \(0\) and \(h\) and \(\eta_{h,k}\) between \(0\) and \(k\)
such that
\[
\frac{R(h,k)}{hk}
=\partial_j\partial_i f
  (x_0+\theta_{h,k}e_i+\eta_{h,k}e_j).
\]
Reversing the order of the two applications gives
\(\widehat\theta_{h,k}\) between \(0\) and \(h\) and
\(\widehat\eta_{h,k}\) between \(0\) and \(k\) for which
\[
\frac{R(h,k)}{hk}
=\partial_i\partial_j f
  (x_0+\widehat\theta_{h,k}e_i+\widehat\eta_{h,k}e_j).
\]
As \((h,k)\to(0,0)\), all the intermediate points converge to
\(x_0\). Continuity of the two mixed derivatives at \(x_0\) proves the
desired equality.
\end{proof}

\begin{corollary}\label{cor:nociones-fundamentales-clase}
 Let $f\colon U\subseteq \mathbb{R}^{n}\longrightarrow \mathbb{R}^{m}$ be of class $C^{2}$ at $x_{0}\in U$.
 Then, for any $i,j\in \{1,\dots,n\}$,
 \[
 \frac{\partial ^{2}f}{\partial x_{i}\partial x_{j}}(x_{0})
 =
 \frac{\partial ^{2}f}{\partial x_{j}\partial x_{i}}(x_{0}).
 \]
\end{corollary}
\begin{proof}
Since \(f\) is of class \(C^2\) at \(x_0\), there is a ball centered at
\(x_0\) on which both mixed derivatives exist and are continuous. Apply
Theorem~\ref{clairaut schwarz parciales conmutan}.
\end{proof}

The preceding notation simplifies higher-order partial derivatives.
We give a few examples of its use.
\begin{example}\label{ej:nociones-fundamentales-consideramos-multiindice-calculamos-paso-paso-usamos}
 \begin{enumerate}
 \item Let $f(x_1,x_2,x_3) = x_1^2 x_2^3 \sin(x_3)$. Consider the multi-index $\alpha = (2,1,1)$. Then
 \[
 D^\alpha f = \frac{\partial^4 f}{\partial x_{1}^2 \partial x_{2} \partial x_{3}}.
 \]
 Differentiating successively, we obtain:
 \[
 \frac{\partial f}{\partial x_{3}} = x_1^2 x_2^3 \cos(x_3), \quad
 \frac{\partial}{\partial x_{2}}\left( \frac{\partial f}{\partial x_{3}} \right) = 3x_1^2 x_2^2 \cos(x_3),
 \]
 \[
 \frac{\partial}{\partial x_{1}}\left( 3x_1^2 x_2^2 \cos(x_3) \right) = 6x_1 x_2^2 \cos(x_3), \quad
 \frac{\partial^2}{\partial x_{1}^2}\left( 3x_1^2 x_2^2 \cos(x_3) \right) = 6 x_2^2 \cos(x_3).
 \]

 \item If $f(x,y) = e^{xy}$ and we use the multi-index $\alpha = (1,2)$, then
 \[
 D^\alpha f = \frac{\partial
 ^3 f}{\partial x_{1} \partial x_{2}^2}=\frac{\partial ^{3}f}{\partial x\partial y^{2}}.
 \]
 We compute:
 \[
 \frac{\partial f}{\partial y} = x e^{xy}, \quad
 \frac{\partial^2 f}{\partial y^2} = x^2 e^{xy}, \quad
 \frac{\partial^3 f}{\partial x \partial y^2} = \frac{\partial}{\partial x}(x^2 e^{xy}) = (2x + x^2 y)e^{xy}.
 \]

 \item Let $f(x_1, x_2) = \ln(1 + x_1^2 + x_2^2)$ and let $\alpha = (1,1)$. Then
 \[
 D^\alpha f = \frac{\partial^2 f}{\partial x_{1} \partial x_{2}}.
 \]
 We compute:
 \[
 \frac{\partial f}{\partial x_{2}} = \frac{2x_2}{1 + x_1^2 + x_2^2}, \quad
 \frac{\partial^2 f}{\partial x_{1} \partial x_{2}} = \frac{\partial}{\partial x_{1}} \left( \frac{2x_2}{1 + x_1^2 + x_2^2} \right) = \frac{-4x_1 x_2}{(1 + x_1^2 + x_2^2)^2}.
 \]

 \item If we take the multi-index $\alpha = (0,\dots,0,k,0,\dots,0)$ with a $k$ in the $i$th position, then
 \[
 D^\alpha f(x) = \frac{\partial^k f}{\partial x_{i}^k}(x),
 \]
 which is simply the derivative of order $k$ with respect to the variable $x^{i}$.
 \end{enumerate}
\end{example}

We use multi-index notation to present the most important properties of higher-order partial derivatives for sufficiently regular functions. Among the properties we will study is the familiar Leibniz rule for products. It is useful first to recall the result for functions from $\mathbb{R}$ to $\mathbb{R}$. Its proof uses properties of binomial coefficients, in particular the following result:
\begin{theorem}[Pascal's triangle]\label{pascal}\index{Pascal's triangle@Pascal's triangle}
For all integers \( n \geq 1 \) and \( 1 \leq k \leq n \), the following identity between binomial coefficients holds:
\[
\binom{n}{k} = \binom{n-1}{k-1} + \binom{n-1}{k}.
\]
This identity reflects the construction of Pascal's triangle: each number is the sum of the two numbers directly above it.
\end{theorem}
\begin{proof}
If \(1\leq k\leq n-1\), then
\[
\binom{n-1}{k-1}+\binom{n-1}{k}
=\frac{k(n-1)!+(n-k)(n-1)!}{k!(n-k)!}
=\frac{n!}{k!(n-k)!}
=\binom nk.
\]
For \(k=n\), using the convention \(\binom{n-1}{n}=0\), the right-hand
side is \(\binom{n-1}{n-1}=1=\binom nn\).
\end{proof}
The Leibniz rule states the following:
\begin{theorem}[generalized Leibniz formula]\label{formula de leibniz orden superior de R en R}\index{generalized Leibniz formula@generalized Leibniz formula}
 Let $n\in\mathbb N$ and $f,g\colon I\subseteq \mathbb{R}\longrightarrow \mathbb{R}$ be $n$-times differentiable functions. Then $fg$ is $n$-times differentiable and
 \[(fg)^{(n)}=\displaystyle\sum_{k=0}^{n}{n\choose k}f^{(k)}g^{(n-k)}.\]
\end{theorem}
\begin{proof}
 We proceed by induction on $n$. The case $n=1$ states that

\[
(fg)' = fg' + f'g,
\]

which is the usual product rule. Suppose that the assertion holds
for a fixed $n\geq1$, and now let $f$ and $g$ be
$(n+1)$-times differentiable functions. The induction hypothesis gives

\[
(fg)^{(n)} = \displaystyle\sum_{k=0}^{n} \binom{n}{k} f^{(k)} g^{(n-k)}.
\]

Then

\[
(fg)^{(n+1)} \overset{\text{induction hypothesis}}{=} \left[ \displaystyle\sum_{k=0}^{n} \binom{n}{k} f^{(k)} g^{(n-k)} \right]'
\]
\[
= \displaystyle\sum_{k=0}^{n} \binom{n}{k} f^{(k)} g^{(n+1-k)}+\displaystyle\sum_{k=0}^{n} \binom{n}{k} f^{(k+1)} g^{(n-k)}
\]
\[
= \displaystyle\sum_{k=0}^{n} \binom{n}{k} f^{(k)} g^{(n+1-k)} + \displaystyle\sum_{k=1}^{n+1} \binom{n}{k-1} f^{(k)} g^{(n+1-k)}
\]
\[
= \binom{n}{0} f^{(0)} g^{(n+1)} + \displaystyle\sum_{k=1}^{n} \left[ \binom{n}{k} + \binom{n}{k-1} \right] f^{(k)} g^{(n+1-k)} + \binom{n}{n} f^{(n+1)} g^{(0)}
\]
\[
\overset{\text{Theorem~\ref{pascal}}}{=} \binom{n+1}{0} f^{(0)} g^{(n+1)} + \displaystyle\sum_{k=1}^{n} \binom{n+1}{k} f^{(k)} g^{(n+1-k)} + \binom{n+1}{n+1} f^{(n+1)} g^{(0)}
\]
\[
= \displaystyle\sum_{k=0}^{n+1} \binom{n+1}{k} f^{(k)} g^{(n+1-k)}.
\]

Thus the assertion holds for \( n+1 \), and the principle of induction shows that it holds for every $n\in \mathbb{N}$.
\end{proof}

For functions of class $C^{|\alpha|}$, multi-index notation gives a formula similar to the preceding one. We conclude this section by presenting this formula together with other properties of higher-order derivatives.
\begin{proposition}[Properties of partial derivatives in multi-index notation]\label{prop: propiedades parciales multiindices incluye leibniz multiindices}\index{properties of partial derivatives in multi-index notation@properties of partial derivatives in multi-index notation}
Let $\alpha, \beta \in \mathbb{N}_0^{n}$,
$f,g\colon U\subseteq\mathbb{R}^n\longrightarrow\mathbb{R}$ be functions of
class $C^{|\alpha|}$, and let
$h\colon U\subseteq\mathbb{R}^n\longrightarrow\mathbb{R}$ be a function of
class $C^{|\alpha+\beta|}$, where
$|\alpha+\beta|=|\alpha|+|\beta|$. Then the following
properties hold:

\begin{enumerate}
 \item Linearity:
 \[
 D^\alpha (\lambda f + \mu g) = \lambda D^\alpha f + \mu D^\alpha g, \quad \text{for any } \lambda,\mu \in \mathbb{R}.
 \]
 \item $D^{\beta}(D^{\alpha}h)=D^{\alpha+\beta}h=D^{\alpha}(D^{\beta}h)$.
 \item Product rule (Leibniz rule):
 \[
 D^\alpha(fg) = \displaystyle\sum_{\beta \leq \alpha} \binom{\alpha}{\beta} D^\beta f \cdot D^{\alpha - \beta} g,
 \]
 where $\displaystyle\binom{\alpha}{\beta} = \prod_{i=1}^n \binom{\alpha_i}{\beta_i}$ is the binomial coefficient and $\beta \leq \alpha$ if and only if $\beta_{j}\leq \alpha_{j}$ for every $j\in\{1,\dots,n\}$.
\end{enumerate}
\end{proposition}

\begin{proof}
 \begin{enumerate}
 \item This follows immediately from the linearity of partial derivatives.
 \item Suppose that $h$ is of class $C^{|\alpha|+|\beta|}$. Repeated
 application of Theorem~\ref{clairaut schwarz parciales conmutan}
 to lower-order derivatives allows us to interchange adjacent partial
 derivatives. Therefore,
 \[
D^{\beta}(D^{\alpha}h)
= \frac{\partial^{|\beta|}}{\partial x_{1}^{\beta_{1}} \cdots \partial x_{n}^{\beta_{n}}}
\left(
 \frac{\partial^{|\alpha|}h}{\partial x_{1}^{\alpha_{1}} \cdots \partial x_{n}^{\alpha_{n}}}
\right)
= \frac{\partial^{|\alpha|+|\beta|}h}{\partial x_{1}^{\beta_{1}} \cdots \partial x_{n}^{\beta_{n}} \partial x_{1}^{\alpha_{1}} \cdots \partial x_{n}^{\alpha_{n}}}
\]
\[
= \frac{\partial^{|\alpha+\beta|}h}{\partial x_{1}^{\alpha_{1}+\beta_{1}} \cdots \partial x_{n}^{\alpha_{n}+\beta_{n}}}
= D^{\alpha+\beta}h.
\]

 which is what we wanted to prove. Interchanging $\alpha$ and $\beta$ in
the same calculation and applying the Clairaut--Schwarz theorem again gives
$D^\alpha(D^\beta h)=D^{\alpha+\beta}h$.
 \item We proceed by induction on $|\alpha|$. For $|\alpha|=0$,
 the sum contains only $\beta=0$ and the identity is $fg=fg$.
 Suppose the formula has been proved for all multi-indices of order
 $r$, and let $|\alpha|=r+1$. Choose $j$ with $\alpha_j>0$ and set
 $\alpha':=\alpha-e_j$. By the induction hypothesis, the usual product
 rule, and the preceding item,
 \[
 \begin{aligned}
 D^\alpha(fg)
 &=D^{e_j}\!\left(
   \displaystyle\sum_{\beta\leq\alpha'}
   \binom{\alpha'}{\beta}D^\beta f\,D^{\alpha'-\beta}g
   \right)\\
 &=\displaystyle\sum_{\beta\leq\alpha'}
   \binom{\alpha'}{\beta}
   \left(
   D^{\beta+e_j}f\,D^{\alpha'-\beta}g
   +D^\beta f\,D^{\alpha'-\beta+e_j}g
   \right).
 \end{aligned}
 \]
 In the first sum we set $\delta=\beta+e_j$ and in the second
 $\delta=\beta$. With the convention that a multi-index
 coefficient is zero when any component of the lower index is
 negative or exceeds the corresponding component of the upper index, the coefficient of
 $D^\delta f\,D^{\alpha-\delta}g$ is
 \[
 \binom{\alpha'}{\delta-e_j}+\binom{\alpha'}{\delta}.
 \]
 The factors corresponding to coordinates other than $j$
 agree, and in coordinate $j$ Pascal's identity gives
 \[
 \binom{\alpha_j-1}{\delta_j-1}
 +\binom{\alpha_j-1}{\delta_j}
 =\binom{\alpha_j}{\delta_j}.
 \]
 Therefore,
 \[
 D^\alpha(fg)
 =\displaystyle\sum_{\delta\leq\alpha}
 \binom{\alpha}{\delta}D^\delta f\,D^{\alpha-\delta}g,
 \]
 which completes the induction.
 \end{enumerate}
\end{proof}

 Motivated by the multinomial theorem, we can state a version of the Leibniz rule for several factors:
\begin{proposition}[Leibniz rule for several factors]
\label{prop:leibniz-multiindices-varios-factores}\index{Leibniz rule for several factors@Leibniz rule for several factors}
Let $\alpha\in\mathbb{N}_{0}^{n}$ and let $m\in\mathbb{N}$ and
$f_{1},\dots,f_{m}\colon U\subseteq\mathbb{R}^{n}\longrightarrow\mathbb{R}$ be functions of class $C^{|\alpha|}$.
Then
\[
D^{\alpha}\!\Bigl(\prod_{r=1}^{m} f_{r}\Bigr)
=
\displaystyle\sum_{\substack{\beta_{1},\dots,\beta_{m}\in\mathbb{N}_{0}^{n}\\
\beta_{1}+\cdots+\beta_{m}=\alpha}}
\frac{\alpha!}{\beta_{1}!\cdots \beta_{m}!}\,
\prod_{r=1}^{m} D^{\beta_{r}} f_{r}.
\]
\end{proposition}
\begin{proof}
We proceed by induction on $m$. For $m=2$, the formula is precisely the
Leibniz rule in Proposition~\ref{prop: propiedades parciales multiindices incluye leibniz multiindices}.

Suppose the assertion holds for $m$ factors and consider
$f_{1},\dots,f_{m+1}$. Write
$F:=\displaystyle\prod_{r=1}^{m} f_{r}$, so that
$\displaystyle\prod_{r=1}^{m+1} f_{r}=F\,f_{m+1}$.
Applying the Leibniz rule for two factors,
\[
D^{\alpha}(F f_{m+1})
=
\displaystyle\sum_{\gamma\leq \alpha}
\binom{\alpha}{\gamma}
D^{\gamma}F\,
D^{\alpha-\gamma}f_{m+1}.
\]

By the induction hypothesis,
\[
D^{\gamma}F
=
\displaystyle\sum_{\substack{\beta_{1},\dots,\beta_{m}\in\mathbb{N}_{0}^{n}\\
\beta_{1}+\cdots+\beta_{m}=\gamma}}
\frac{\gamma!}{\beta_{1}!\cdots\beta_{m}!}
\prod_{r=1}^{m} D^{\beta_{r}} f_{r}.
\]

Substituting and using $\binom{\alpha}{\gamma}\gamma!
=
\displaystyle\frac{\alpha!}{(\alpha-\gamma)!}$, we obtain
\[
D^{\alpha}\left(\prod_{r=1}^{m+1} f_{r}\right)
=
\displaystyle\sum_{\gamma\leq\alpha}
\displaystyle\sum_{\substack{\beta_{1},\dots,\beta_{m}\in\mathbb{N}_{0}^{n}\\
\beta_{1}+\cdots+\beta_{m}=\gamma}}
\frac{\alpha!}{\beta_{1}!\cdots\beta_{m}!(\alpha-\gamma)!}
\left(\prod_{r=1}^{m} D^{\beta_{r}} f_{r}\right)
D^{\alpha-\gamma}f_{m+1}.
\]

Now define
\[
\beta_{m+1}:=\alpha-\gamma.
\]
Since $\gamma\leq\alpha$, we have $\beta_{m+1}\in\mathbb{N}_{0}^{n}$ and,
conversely, given any family
$\beta_{1},\dots,\beta_{m+1}\in\mathbb{N}_{0}^{n}$ such that
\[
\beta_{1}+\cdots+\beta_{m+1}=\alpha,
\]
if we set
\[
\gamma:=\beta_{1}+\cdots+\beta_{m},
\]
then $\gamma\leq\alpha$ and
$\beta_{m+1}=\alpha-\gamma$.
This correspondence is a bijection between the pairs
\[
(\gamma,\beta_{1},\dots,\beta_{m})
\quad\text{with}\quad
\beta_{1}+\cdots+\beta_{m}=\gamma\leq\alpha
\]
and the families
\[
(\beta_{1},\dots,\beta_{m+1})
\quad\text{with}\quad
\beta_{1}+\cdots+\beta_{m+1}=\alpha.
\]

Reindexing the sum by this bijection gives
\[
D^{\alpha}\left(\prod_{r=1}^{m+1} f_{r}\right)
=
\displaystyle\sum_{\substack{\beta_{1},\dots,\beta_{m+1}\in\mathbb{N}_{0}^{n}\\
\beta_{1}+\cdots+\beta_{m+1}=\alpha}}
\frac{\alpha!}{\beta_{1}!\cdots\beta_{m+1}!}
\prod_{r=1}^{m+1} D^{\beta_{r}} f_{r},
\]
which completes the induction.

\end{proof}
The following multi-index form of the chain rule gives explicit control
over the partial derivatives of a composition.
\begin{theorem}[Faà di Bruno formula for multi-indices]\label{faa di bruno multivariable}\index{Faa di Bruno formula for multi-indices@Faà di Bruno formula for multi-indices}
Let \(f\colon U\subseteq\mathbb R^n\longrightarrow\mathbb R^m\) and
\(g\colon V\subseteq\mathbb R^m\longrightarrow\mathbb R\), with
\(f(U)\subseteq V\), both be of class \(C^k\). Write
\(f=(f_1,\dots,f_m)\). Then \(g\circ f\) is of class \(C^k\) and,
for \(x_0\in U\) and
\(\alpha\in\mathbb N_0^n\) with \(1\leq|\alpha|\leq k\),
 \[
D^{\alpha}(g\circ f)(x_0) =
\displaystyle\sum_{\substack{\beta\in \mathbb{N}_{0}^{m},
\text{ }1\leq|\beta|\leq |\alpha|\\
\gamma=(\gamma_{1},\dots,\gamma_{|\beta|}),\text{ }\gamma_{i}\in \mathbb{N}_{0}^{n},\text{ }|\gamma_{i}|>0,\text{ }\sum_{i=1}^{|\beta|}\gamma_{i}=\alpha \\
l=(l_{1},\dots,l_{|\beta|}),\text{ }l_{i}\in\{1,\dots,m\}\\
\#\{i\mid l_i=j\}=\beta_j,\text{ }1\leq j\leq m
}}c_{\alpha,\beta,\gamma,l}D^{\beta}g(f(x_{0}))\prod_{i=1}^{|\beta|}D^{\gamma_{i}}f_{l_{i}}(x_{0}).
\]
The coefficients \(c_{\alpha,\beta,\gamma,l}\) are nonnegative
combinatorial constants depending only on the indicated
indices.
\end{theorem}
\begin{proof}
We proceed by induction on \(q:=|\alpha|\). If \(q=1\), there exists
\(r\in\{1,\dots,n\}\) with \(\alpha=e_r\), and the chain rule gives
\[
D^{e_r}(g\circ f)(x_0)
=\displaystyle\sum_{j=1}^m(D^{e_j}g)(f(x_0))D^{e_r}f_j(x_0).
\]
This is the formula in the statement with \(\beta=e_j\),
\(\gamma=(e_r)\), \(l=(j)\), and coefficient \(1\).

Suppose the assertion holds for multi-indices of length \(q\), and let
\(\alpha'\) have length \(q+1\). Choose \(r\) with \(\alpha'_r>0\) and
set \(\alpha=\alpha'-e_r\). Each term in the inductive formula has
the form
\[
c\,(D^\beta g)(f(x))
\prod_{a=1}^{|\beta|}D^{\gamma_a}f_{l_a}(x),
\qquad
\displaystyle\sum_{a=1}^{|\beta|}\gamma_a=\alpha.
\]
Differentiating with respect to \(x^r\), the derivative of the first factor yields
\[
c\displaystyle\sum_{j=1}^m(D^{\beta+e_j}g)(f(x))D^{e_r}f_j(x)
\prod_{a=1}^{|\beta|}D^{\gamma_a}f_{l_a}(x).
\]
This adds \(e_r\) to the family \((\gamma_a)\), adds \(j\) to
\((l_a)\), and replaces \(\beta\) by \(\beta+e_j\). If the derivative acts
on the factor with index \(a\), it yields
\[
c\,(D^\beta g)(f(x))D^{\gamma_a+e_r}f_{l_a}(x)
\prod_{\substack{b=1\\ b\ne a}}^{|\beta|}
D^{\gamma_b}f_{l_b}(x).
\]
Here \(\gamma_a\) is replaced by \(\gamma_a+e_r\). In both cases, the
new multi-indices sum to \(\alpha'\), and all the
restrictions in the statement are preserved. Like terms are collected;
their coefficients are added and remain nonnegative combinatorial
constants. This completes the induction. The same construction shows, order
by order, that the resulting derivatives are continuous, so
\(g\circ f\) is of class \(C^k\).
\end{proof}

For vector-valued compositions it is useful to retain the directions
of differentiation. If $f\colon U\subseteq\mathbb R^n\to\mathbb R^m$
is of class $C^k$, write
\[
 D^qf(x)[v_1,\ldots,v_q]
 :=D_{v_q,\ldots,v_1}f(x),\qquad 1\leq q\leq k.
\]
This expression is linear in each direction and symmetric. Indeed,
if $v_a=(v_a^1,\ldots,v_a^n)$, the properties of partial derivatives give
\[
 D^qf(x)[v_1,\ldots,v_q]
 =\sum_{i_1,\ldots,i_q=1}^{n}
 v_1^{i_1}\cdots v_q^{i_q}
 \partial_{i_q}\cdots\partial_{i_1}f(x).
\]

\begin{theorem}[Faà di Bruno's formula using partitions]
\label{teo:faa-di-bruno-particiones}
\index{Faa di Bruno formula@Faà di Bruno's formula!using partitions}
Let $k\geq1$ be an integer, let $U\subseteq\mathbb R^n$ and
$V\subseteq\mathbb R^m$ be open, and let $f\in C^k(U,\mathbb R^m)$
and $g\in C^k(V,\mathbb R^d)$ satisfy $f(U)\subseteq V$.
For $1\leq q\leq k$, let $\Pi_q$ be the set of partitions of
$\{1,\ldots,q\}$ into nonempty, pairwise disjoint blocks whose
union is $\{1,\ldots,q\}$.
Write each partition as $\pi=\{B_1,\ldots,B_a\}$, where
$a=|\pi|$ and $\min B_1<\cdots<\min B_a$.
For $v=(v_1,\ldots,v_q)$ and $B=\{i_1<\cdots<i_b\}$, define
\[
 D_Bf(x)[v]:=D^bf(x)[v_{i_1},\ldots,v_{i_b}].
\]
Then, for $x\in U$ and $v_1,\ldots,v_q\in\mathbb R^n$,
\begin{equation}
\label{eq:faa-di-bruno-particiones}
\begin{split}
 D^q(g\circ f)(x)[v_1,\ldots,v_q]
 =\sum_{\substack{\pi\in\Pi_q\\\pi=\{B_1,\ldots,B_a\}}}
 D^ag(f(x))\bigl[D_{B_1}f(x)[v],\ldots,D_{B_a}f(x)[v]\bigr].
\end{split}
\end{equation}
In each summand, $a$ is the number of blocks of the corresponding
partition. The brackets contain a list of $a$ vectors in
$\mathbb R^m$, on which the $a$-linear map $D^ag(f(x))$ is
evaluated. Symmetry of this map allows the blocks to be reordered
without changing the value of the summand.

In particular, if $I\subseteq\mathbb R$ is an open interval and
$f\in C^k(I,\mathbb R^m)$ satisfies $f(I)\subseteq V$, then, for
$t\in I$ and $1\leq q\leq k$,
\begin{equation}
\label{eq:faa-di-bruno-curva-coeficientes}
\begin{split}
 (g\circ f)^{(q)}(t)
 =\sum_{j=1}^{q}
 \sum_{\substack{k_1,\ldots,k_j\geq1\\k_1+\cdots+k_j=q}}
 \frac{q!}{j!\,k_1!\cdots k_j!}\,
 D^jg(f(t))\bigl[f^{(k_1)}(t),\ldots,f^{(k_j)}(t)\bigr],
\end{split}
\end{equation}
where the indices $k_1,\ldots,k_j$ are integers. If both maps
are smooth, the identities hold for every integer $q\geq1$.
\end{theorem}

\begin{proof}
The chain rule, applied successively to the components, shows
that $g\circ f$ is of class $C^k$. We prove the first identity
by induction on $q$. For $q=1$, the only partition is
$\{\{1\}\}$, and the identity reduces to
\[
 D(g\circ f)(x)[v_1]=Dg(f(x))[Df(x)[v_1]].
\]

Suppose the identity holds for an integer $q<k$.
Fix constant directions $v_1,\ldots,v_{q+1}$ and a partition
$\pi=\{B_1,\ldots,B_a\}\in\Pi_q$.
For $1\leq\ell\leq a$, put
\[
 W_\ell(x):=D_{B_\ell}f(x)[v],\qquad
 W_\ell^+(x):=D^{|B_\ell|+1}f(x)
       [v_{i_1},\ldots,v_{i_{|B_\ell|}},v_{q+1}],
 \quad B_\ell=\{i_1<\cdots<i_{|B_\ell|}\}.
\]
Since $|B_\ell|\leq q<k$, the map $W_\ell$ is of class $C^1$,
and its derivative in direction $v_{q+1}$ is $W_\ell^+$.
The map $x\mapsto D^ag(f(x))$ is also of class $C^1$, since
$a\leq q<k$. These derivatives are taken in $U$: for each
$x\in U$, there is $\delta>0$ such that $x+t v_{q+1}\in U$
whenever $|t|<\delta$.

Consider the summand
$T_\pi(x):=D^ag(f(x))[W_1(x),\ldots,W_a(x)]$.
Writing this evaluation in components gives a finite sum of
products of one component of $D^ag(f(x))$ and one component of
each $W_\ell(x)$. The product rule for these products and the
chain rule for the first factor give
\[
\begin{split}
 D_{v_{q+1}}T_\pi(x)
 ={}&D^{a+1}g(f(x))
       [W_1(x),\ldots,W_a(x),Df(x)[v_{q+1}]]\\
 &+\sum_{\ell=1}^{a}D^ag(f(x))
       [W_1(x),\ldots,W_{\ell-1}(x),W_\ell^+(x),
        W_{\ell+1}(x),\ldots,W_a(x)].
\end{split}
\]
The first term corresponds to adjoining the block $\{q+1\}$ to
$\pi$. The term with index $\ell$ in the sum corresponds to
replacing $B_\ell$ by $B_\ell\cup\{q+1\}$ and leaving the
other blocks unchanged. If their order changes, symmetry of the
derivative of $g$ restores the order fixed in the statement.

We check that no partition is omitted or counted twice.
Let $\rho\in\Pi_{q+1}$, and let $C$ be its unique block
containing $q+1$. If $C=\{q+1\}$, removing $C$ leaves a unique
partition $\pi\in\Pi_q$; hence $\rho$ comes from the first term
associated with this $\pi$. If $|C|\geq2$, replacing $C$ by
$C\setminus\{q+1\}$ leaves a unique partition $\pi\in\Pi_q$
and a unique block of $\pi$ to which $q+1$ must be added. Thus
$\rho$ comes from exactly one term of the sum over $\ell$.
The two cases are disjoint. Summing the preceding identity over
$\Pi_q$ therefore gives precisely
\eqref{eq:faa-di-bruno-particiones} for $q+1$. This completes
the induction.

For the curve formula, take $v_1=\cdots=v_q=1$.
A block of size $k_\ell$ then contributes the vector
$f^{(k_\ell)}(t)$. Fix the number $j$ of blocks. We may sum over
all partitions with ordered blocks and divide by $j!$: the
blocks are distinct subsets, so each partition has exactly
$j!$ orderings, even when some blocks have equal sizes. The
evaluation of $D^jg(f(t))$ has the same value for all these
orderings.

For ordered sizes $k_1,\ldots,k_j\geq1$ with sum $q$, choose
the $k_1$ positions of the first block, then the $k_2$ positions
of the second among the remaining positions, and continue.
The number of choices is
\[
 \binom{q}{k_1}\binom{q-k_1}{k_2}\cdots
 \binom{q-k_1-\cdots-k_{j-1}}{k_j}
 =\frac{q!}{k_1!\cdots k_j!}.
\]
All these choices produce the same summand. Multiplying by
this number, dividing by $j!$, and summing over all sizes and
all values $1\leq j\leq q$ yields
\eqref{eq:faa-di-bruno-curva-coeficientes}.
\end{proof}

We make the relation to Theorem~\ref{faa di bruno multivariable}
explicit. Let $\alpha\in\mathbb N_0^n$, with $q=|\alpha|\geq1$,
and choose a list $v_1,\ldots,v_q$ of standard basis vectors in
which $e_s$ occurs $\alpha_s$ times for each $1\leq s\leq n$.
For each block $B_\ell$, define $\gamma_\ell\in\mathbb N_0^n$
by counting how many of its positions correspond to each $e_s$.
Then $|\gamma_\ell|=|B_\ell|\geq1$,
$\displaystyle\sum_{\ell=1}^{a}\gamma_\ell=\alpha$, and
$D_{B_\ell}f(x)[v]=D^{\gamma_\ell}f(x)$.
If $g$ is scalar-valued, expanding the multilinear evaluation gives
\[D^ag(f(x))
   [D^{\gamma_1}f(x),\ldots,D^{\gamma_a}f(x)]=\sum_{l_1,\ldots,l_a=1}^{m}
 (\partial_{l_a}\cdots\partial_{l_1}g)(f(x))\prod_{\ell=1}^{a}D^{\gamma_\ell}f_{l_\ell}(x).
\]
For each list $(l_1,\ldots,l_a)$, let $\beta_t$ be the number
of indices $\ell$ such that $l_\ell=t$, for $1\leq t\leq m$.
Then $|\beta|=a$, and the first factor is
$(D^\beta g)(f(x))$. Grouping pairs consisting of a
partition and a component list that produce the same indices
$(\beta,\gamma,l)$ gives nonnegative integer coefficients
$c_{\alpha,\beta,\gamma,l}$; a coefficient is zero if no pair
produces those indices. This recovers the form of the preceding
theorem. Coefficients can be distributed in different ways among
lists that differ only by a permutation of factors; the sum and
the derivative it represents are unchanged.

\subsection{The Banach fixed point theorem and the inverse function theorem}

\begin{theorem}[Banach fixed point theorem]
\label{teo:principio-contraccion-banach}
\index{Banach fixed point theorem}\index{contraction principle}
Let $(M,d)$ be a nonempty complete metric space, and let $T\colon M\longrightarrow M$ be a contraction; that is, suppose there exists $q\in[0,1)$ such that $d(Tx,Ty)\leq qd(x,y)$ for all $x,y\in M$. Then $T$ has a unique fixed point $x_*$. Moreover, for any $x_0\in M$, the sequence defined by $x_{n+1}:=T(x_n)$ converges to $x_*$ and
\[
d(x_n,x_*)\leq\frac{q^n}{1-q}d(x_1,x_0)
\]
for every integer $n\geq1$.
\end{theorem}

\begin{proof}
If $q=0$, the map $T$ is constant and its value is the unique fixed point;
the sequence reaches that point after its first iteration.
Assume $0<q<1$. Fix $x_0\in M$. The contraction inequality gives $d(x_2,x_1)\leq qd(x_1,x_0)$. If $d(x_{n+1},x_n)\leq q^nd(x_1,x_0)$ for some $n\in\mathbb N$, applying $T$ to the two points gives $d(x_{n+2},x_{n+1})\leq qd(x_{n+1},x_n)\leq q^{n+1}d(x_1,x_0)$. By induction, the estimate holds for every $n\in\mathbb N$. If $m>n$, the triangle inequality gives
\[
d(x_m,x_n)
\leq
\displaystyle\sum_{j=n}^{m-1}d(x_{j+1},x_j)
\leq
\displaystyle\sum_{j=n}^{m-1}q^jd(x_1,x_0)
\leq
\frac{q^n}{1-q}d(x_1,x_0).
\]
Thus $(x_n)$ is a Cauchy sequence. Completeness of $M$ guarantees existence of $x_*\in M$ such that $x_n\longrightarrow x_*$. Since every contraction is continuous,
$\displaystyle T(x_*)=\lim_{n\to\infty}T(x_n)=\lim_{n\to\infty}x_{n+1}=x_*$. If $y$ is another fixed point, then $d(x_*,y)=d(Tx_*,Ty)\leq qd(x_*,y)$; since $q<1$, it follows that $d(x_*,y)=0$. Finally, letting $m\to\infty$ in the preceding estimate gives the inequality in the statement.
\end{proof}

\begin{theorem}[Inverse function theorem]
\label{lem:funcion-inversa-euclidea-previa-variedades}
\index{inverse function theorem!in Euclidean spaces}
Let $U\subseteq\mathbb R^n$ be open, let $f\in C^j(U,\mathbb R^n)$,
where $j\in\mathbb N$, $j\geq1$, or $j=\infty$, and let $a\in U$.
If $Df(a)$ is invertible, there exist connected open sets $V\subseteq U$
and $W\subseteq\mathbb R^n$, with $a\in V$ and $f(a)\in W$, such that
$f\restriction_V\colon V\longrightarrow W$ is a diffeomorphism of class
$C^j$: it is bijective, and both it and its inverse $g$ are of class $C^j$.
Moreover,
\[
Dg(y)=Df(g(y))^{-1},\qquad y\in W.
\]
\end{theorem}

\begin{proof}
If $n=0$, the spaces consist of a single point and the result is immediate.
Assume $n\geq1$. First consider the case $a=0$, $f(0)=0$, and $Df(0)=I$.
Continuity of $Df$ gives $\rho>0$ such that
$\overline B_{\mathrm{euc}}(0,\rho)\subseteq U$ and
\begin{equation}
\label{eq:contraccion-prueba-inversa-euclidea}
\|Df(x)-I\|_{\operatorname{op}}\leq\frac12,
\qquad x\in\overline B_{\mathrm{euc}}(0,\rho).
\end{equation}
Set $R(x):=f(x)-x$. For $x,z$ in this closed ball, the fundamental theorem
of calculus, applied componentwise along the segment from $z$ to $x$, gives
\[
R(x)-R(z)=\int_0^1\bigl(Df(z+t(x-z))-I\bigr)(x-z)\,dt.
\]
Consequently,
\[
\|R(x)-R(z)\|\leq\frac12\|x-z\|,
\qquad
\|x-z\|\leq2\|f(x)-f(z)\|.
\]
The second inequality follows by writing
$x-z=f(x)-f(z)-(R(x)-R(z))$. In particular, $f$ is injective on the ball.

Let $W:=B_{\mathrm{euc}}(0,\rho/2)$ and
$V:=f^{-1}(W)\cap B_{\mathrm{euc}}(0,\rho)$. For $y\in W$, define
\[
T_y(x):=x+y-f(x)=y-R(x),
\qquad x\in\overline B_{\mathrm{euc}}(0,\rho).
\]
Since $R(0)=0$, we have
\[
\|T_y(x)\|\leq\|y\|+\frac12\|x\|<\rho,
\qquad
\|T_y(x)-T_y(z)\|\leq\frac12\|x-z\|.
\]
The closed ball is a complete metric space.
Theorem~\ref{teo:principio-contraccion-banach} provides a unique fixed
point $g(y)$ of $T_y$ in this ball. The first inequality shows that
$g(y)\in B_{\mathrm{euc}}(0,\rho)$, and the fixed point equation is
equivalent to $f(g(y))=y$. Thus $f\restriction_V$ is bijective onto $W$,
with inverse $g$. Moreover,
\begin{equation}
\label{eq:lipschitz-inversa-euclidea-local}
\|g(y)-g(z)\|\leq2\|y-z\|,\qquad y,z\in W.
\end{equation}
Hence $g$ is continuous. The set $V$ is open and, since $V=g(W)$ and
$W$ is connected, it is also connected.

For $x\in V$, the estimate
$\|Df(x)v-v\|\leq\tfrac12\|v\|$ implies that $Df(x)v=0$ only if $v=0$.
Since $Df(x)$ is an endomorphism of $\mathbb R^n$, it is invertible.
Fix $y\in W$, set $x:=g(y)$, and, for $h\ne0$ with $y+h\in W$, define
$u:=g(y+h)-g(y)$. Then $u\ne0$, $\|u\|\leq2\|h\|$, and
\[
g(y+h)-g(y)-Df(x)^{-1}h
=-Df(x)^{-1}\bigl(f(x+u)-f(x)-Df(x)u\bigr).
\]
Consequently,
\[
\frac{\|g(y+h)-g(y)-Df(x)^{-1}h\|}{\|h\|}
\leq2\|Df(x)^{-1}\|_{\operatorname{op}}
\frac{\|f(x+u)-f(x)-Df(x)u\|}{\|u\|}
\longrightarrow0
\]
as $h\to0$, by differentiability of $f$ at $x$. This proves
$Dg(y)=Df(g(y))^{-1}$.

The inversion map $\operatorname{Inv}(A):=A^{-1}$ is smooth on
$\operatorname{GL}(n,\mathbb R)$:
the formula $A^{-1}=\operatorname{adj}(A)/\det A$ expresses each entry
as a quotient of polynomials with a nonzero denominator. Continuity of
$g$ and $Df$ therefore shows that $Dg$ is continuous and $g\in C^1$.
If $f\in C^j$ and $g\in C^\ell$ with $1\leq\ell<j$, the chain rule
implies that $\operatorname{Inv}\circ Df\circ g$ is of class $C^\ell$.
Since this map is $Dg$, it follows that $g\in C^{\ell+1}$. Induction
proves $C^j$ regularity for every finite $j$; if $j=\infty$, the same
argument applies at each finite order.

For the general case, let $A:=Df(a)$ and define
\[
F(x):=A^{-1}\bigl(f(a+x)-f(a)\bigr),\qquad x\in U-a.
\]
Then $F(0)=0$ and $DF(0)=I$. The preceding case provides connected open
sets $V'$, $W'$ and an inverse $G\colon W'\longrightarrow V'$ of class
$C^j$. Taking
\[
V:=a+V',\qquad W:=f(a)+A(W'),\qquad
g(y):=a+G\bigl(A^{-1}(y-f(a))\bigr),
\]
we obtain the inverse of $f\restriction_V$ with the required regularity.
The chain rule applied to $f\circ g=\operatorname{Id}_W$ gives the stated
formula for $Dg$.
\end{proof}

\section{Smooth manifolds}

A smooth manifold formalizes the idea of a space that may be globally curved or topologically complicated, yet admits Euclidean coordinates around each point. The sphere and the torus are familiar examples: no single chart describes either one completely, but every sufficiently small region can be compared with an open subset of the plane.

The dimension of these spaces is intrinsic. Although the usual sphere is drawn inside $\mathbb{R}^{3}$, only two local coordinates are needed to move along it. Requiring a space to be locally homeomorphic to an open subset of $\mathbb{R}^{n}$ captures this topological observation; requiring smooth compatibility between charts then makes it possible to transport calculus.
\begin{figure}[htbp]
    \centering
    \includegraphics[height=4cm]{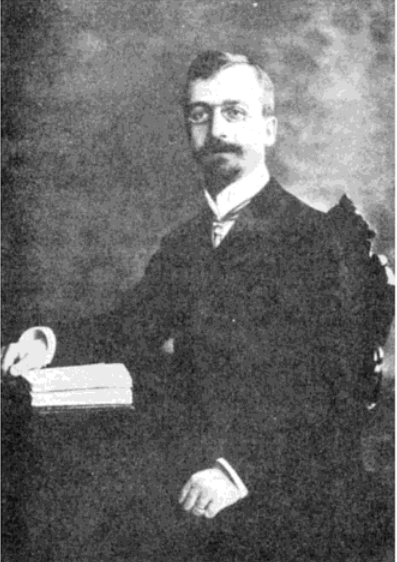}

    \smallskip
    \footnotesize Élie Cartan (1869--1951)

    \caption[Élie Cartan]{Élie Cartan was a central figure in the development
    of modern differential geometry. His work on differential forms,
    moving frames, connections, and curvature made decisive contributions to the
    intrinsic study of manifolds and established part of the language that now
    underlies Riemannian geometry and global analysis.}
    \label{fig:elie-cartan}
\end{figure}
This concept is fundamental in nonlinear analysis, since many partial differential equations and variational problems are naturally formulated on spaces that are not open subsets of $\mathbb{R}^n$. For example:
\medskip
\begin{itemize}
 \item Models in mathematical physics whose domain is a curved surface or a space with variable curvature (waves on the sphere, diffusion on the torus).
 \item Problems in differential geometry involving the boundary of a manifold, such as trace theorems for Sobolev functions defined on the boundary.
 \item Elliptic and parabolic equations on compact spaces with or without boundary, where the Riemannian metric directly influences the differential operators.
\end{itemize}

The notion of a smooth manifold allows us to extend differential calculus to these settings. Once manifolds are equipped with a Riemannian metric, we will have the tools to define gradients, divergences, Laplacians, and elliptic operators intrinsically. This will allow us to develop the study of Sobolev spaces on manifolds, state and prove trace theorems, and analyze partial differential equations in a general geometric setting.

We will present the modern definition of a \emph{differentiable manifold}, first introducing several preliminary concepts.

\begin{definition}\label{def:nociones-fundamentales-variedad-topologica-de-dimension}\index{topological manifold of dimension@topological manifold of dimension}
 Let $M$ be a topological space. We say that $M$ is a \textbf{topological manifold without boundary of dimension $n$} if:
 \begin{enumerate}[label=(\alph*)]
 \item $M$ is Hausdorff.
 \item $M$ is second countable.
 \item For every $p\in M$, there exist an open set $U\subseteq M$ and a homeomorphism $\phi\colon U\longrightarrow \phi(U)$, where $\phi(U)\subseteq \mathbb{R}^{n}$ is open ($M$ is locally homeomorphic to $\mathbb{R}^{n}$), such that $p\in U$.
 \end{enumerate}
\end{definition}

To become familiar with the property of \emph{being locally homeomorphic to $\mathbb{R}^{n}$}, it is useful to consider the following picture:

\begin{figure}[h]

 \includegraphics[scale=0.2]{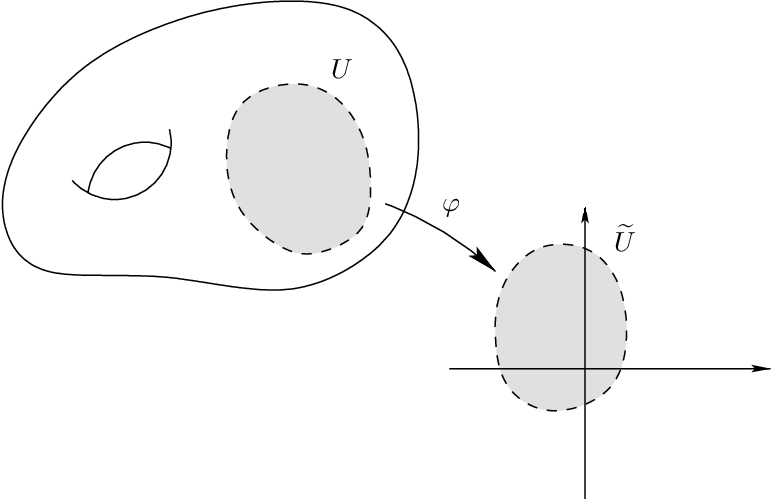}
\centering
\label{fig:variedad-localmente-euclidiana}
\end{figure}
\begin{remark}\label{obs:nociones-fundamentales-carta-coordenada}
 To simplify notation, the open set $U$ and the homeomorphism \[\phi\colon U\longrightarrow \phi(U)\subseteq\mathbb{R}^{n}\] are usually considered together as an ordered pair $(U,\phi)$, called a \emph{coordinate chart}, or simply a \emph{chart}.
\end{remark}
The topological structure captures the locally Euclidean character of the space, but does not yet allow differentiation. To transport calculus between charts, we must require the changes of coordinates to be smooth. A compatible atlas thus assembles the local information into a coherent global structure.

\begin{definition}\label{def:nociones-fundamentales-variedad}\index{atlas smooth structure and chart compatibility@atlas, smooth structure, and chart compatibility}
 Let $M$ be a topological manifold without boundary of dimension $n$.
 \begin{enumerate}[label=(\alph*)]
 \item We say that two charts $(U,\phi)$ and $(V,\psi)$ are smoothly compatible if $U\cap V= \varnothing$ or if the transition map \[\psi\circ \phi^{-1}\colon \phi(U\cap V)\longrightarrow \psi(U\cap V)\] is a diffeomorphism of class $C^{\infty}$ (henceforth called smooth).
 \item An atlas for $M$ is a collection of charts $\mathcal{A}$ such that $\displaystyle\bigcup_{(U,\phi)\in \mathcal{A}}U=M$.
 \item A smooth atlas is an atlas $\mathcal{A}$ such that any two charts in $\mathcal{A}$ are smoothly compatible.
 \item A smooth structure on $M$ is a smooth atlas $\mathcal{A}$ maximal with respect to $\subseteq$ among the smooth atlases of $M$.
 \end{enumerate}
\end{definition}
We can now define a \emph{differentiable manifold (or smooth manifold)}:
\begin{definition}\label{def:nociones-fundamentales-variedad-2}\index{smooth manifold}
 A smooth manifold without boundary of dimension $n$ is a topological manifold $M$ of dimension $n$ equipped with a smooth structure $\mathcal{A}$. We write $(M,\mathcal{A})$ when we wish to refer explicitly to the smooth structure.
\end{definition}
\begin{remark}\label{obs:nociones-fundamentales-tambien-definir-variedades-clase-aquellas-cuales}
 We can also define manifolds without boundary of class $C^{k}$ by requiring the transition maps \[\psi\circ\phi^{-1}\colon \phi(U\cap V)\longrightarrow \psi(U\cap V)\] to be diffeomorphisms of class $C^{k}$, that is, functions of class $C^{k}$ whose inverses are of class $C^{k}$.
\end{remark}

\subsection{Examples of manifolds and construction of their charts}

\begin{example}[Euclidean open sets]
\label{ej:variedades-abiertos-euclidianos}
Let $U\subseteq\mathbb R^n$ be open, with its subspace topology. It is
Hausdorff, and intersections with $U$ of balls with rational centers
and radii form a countable basis. The identity
$\operatorname{id}_U\colon U\to U$ is a global chart. Its restrictions
to open subsets of $U$ are compatible, since the coordinate changes
are restrictions of the identity. The maximal atlas they generate
makes $U$ a smooth manifold of dimension $n$. In particular, both
$\mathbb R^n$ and the open ball are manifolds without boundary: the
boundary of a subset in an ambient space and the boundary of a manifold
are different concepts.
\end{example}

\begin{example}[The sphere through hemisphere charts]
\label{ej:variedades-esfera-cartas-hemisferios}
For $n\geq1$, consider
\[
 S^n:=\left\{x\in\mathbb R^{n+1}\middle|
                 \sum_{j=1}^{n+1}(x^j)^2=1\right\}.
\]
It is Hausdorff and second countable with its subspace topology; it is
compact because it is closed and bounded. For
$i\in\{1,\ldots,n+1\}$ and $\varepsilon\in\{-1,1\}$, set
\[
 U_{i,\varepsilon}:=\{x\in S^n\mid \varepsilon x^i>0\},
 \qquad
 \phi_{i,\varepsilon}(x):=(x^1,\ldots,\widehat{x^i},\ldots,x^{n+1}).
\]
The hat means that this coordinate is deleted. The image is
$B^n:=\{y\in\mathbb R^n\mid |y|<1\}$. The inverse inserts, in
position $i$, the coordinate
\[
 x^i=\varepsilon\sqrt{1-|y|^2}.
\]
The projection and this inverse are continuous, so each pair is a chart.
The $2(n+1)$ domains cover the sphere, since every unit vector has a
nonzero component. On an overlap, the transition
$\phi_{j,\delta}\circ\phi_{i,\varepsilon}^{-1}$ is obtained by inserting
the displayed square root and deleting coordinate $j$. Its components
are coordinates of $y$ or the function
$\varepsilon\sqrt{1-|y|^2}$, which are smooth on $B^n$. The inverse has
the same description with the indices exchanged. Thus this atlas defines
a compact smooth manifold without boundary of dimension $n$.
The equation defining $S^n$ in the ambient space is not a chart:
the charts are the explicit maps just constructed.
\end{example}

\begin{example}[Real projective space]
\label{ej:variedades-proyectivo-real-cartas}
For $n\geq1$, let $\mathbb{RP}^n$ be the set of vector lines in
$\mathbb R^{n+1}$, with the quotient topology of
$q\colon S^n\to\mathbb{RP}^n$, $q(x)=\mathbb Rx$. It is compact as the
continuous image of $S^n$. To verify the Hausdorff property without
assuming it for the quotient, consider
\[
 P\colon\mathbb{RP}^n\longrightarrow\mathbb R^{(n+1)\times(n+1)},
 \qquad P(\mathbb Rx):=xx^{\mathsf T},\quad |x|=1.
\]
Replacing $x$ by $-x$ does not change the matrix. The equality
$xx^{\mathsf T}=yy^{\mathsf T}$ implies, by applying both matrices to
$x$, that $x=\langle y,x\rangle y$; hence the lines agree.
Thus $P$ is injective. It is continuous by the definition of the quotient
topology, since $P\circ q$ is continuous. A continuous bijection from a
compact space to a Hausdorff space is a homeomorphism onto its image:
the image of every closed subset is compact, hence closed.
Applying this observation to $P$ proves that $\mathbb{RP}^n$ is
Hausdorff.

For $i\in\{1,\ldots,n+1\}$, let $V_i$ be the set of lines admitting
a representative with $x^i\ne0$. It is open, since its preimage under $q$
is open. Define
\[
 \kappa_i(\mathbb Rx):=
 \left(\frac{x^1}{x^i},\ldots,
       \widehat{\frac{x^i}{x^i}},\ldots,
       \frac{x^{n+1}}{x^i}\right)\in\mathbb R^n.
\]
The ratios do not depend on the representative. This map is continuous,
and its inverse sends $y$ to the line of the vector $a(y)$ whose
component $i$ is one and whose other components are those of $y$.
Indeed, the inverse is $q(a(y)/|a(y)|)$, which is continuous. On the
domain of $\kappa_k\circ\kappa_i^{-1}$, we have $a^k(y)\ne0$, and its
components are
\[
 \frac{a^j(y)}{a^k(y)},\qquad j\ne k.
\]
They are smooth wherever the denominator is nonzero; the inverse
transition has the same form. The $n+1$ charts cover the space. Their
countable bases, transported and combined, provide a countable basis
for $\mathbb{RP}^n$. We obtain a compact smooth manifold without boundary
of dimension $n$, constructed from a quotient and charts without
choosing global coordinates on the manifold.
\end{example}

\subsection{Topological properties and manifolds with boundary}

Topological manifolds (and therefore smooth manifolds) have several interesting topological properties:
\begin{proposition}\label{prop:nociones-fundamentales-variedad-topologica-localmente-conexo-trayectorias-conexa}
 Let $M$ be a topological manifold without boundary.
 \begin{enumerate}
 \item $M$ is locally path connected.
 \item $M$ is connected if and only if it is path connected.
 \item The connected components of $M$ coincide with its path components.
 \item $M$ has countably many connected components, each of which is an open subset of $M$ and a connected topological manifold.
 \end{enumerate}

\end{proposition}
\begin{proof}
Let \(p\in M\) and let \(U\) be an open neighborhood of \(p\). Choose a
chart \((V,\phi)\) with \(p\in V\). The open set
\(\phi(U\cap V)\) contains an open Euclidean ball \(B\) centered at
\(\phi(p)\). Then \(\phi^{-1}(B)\) is a neighborhood of \(p\), is
contained in \(U\), and is path connected because \(B\) is convex.
This proves (1).

Now fix \(p\in M\) and let \(C_p\) be the set of points that can
be joined to \(p\) by a path. By (1), \(C_p\) is open.
The same argument at each point of \(M\setminus C_p\) shows that
its complement is also open. If \(M\) is connected, it follows that
\(C_p=M\). Since every path connected space is connected, this
proves (2). Applying (2) to each connected component gives (3), and
local path connectedness also shows that every component is
open.

Finally, let \(\mathcal B\) be a countable basis for \(M\). Each connected
component \(C\) is open and nonempty, and therefore contains some nonempty basis
 element \(B_C\in\mathcal B\). Distinct components contain distinct basis
elements, since they are disjoint; hence the family of components is
countable. Restricting the charts of \(M\) makes each
component a connected topological manifold. This proves (4).
\end{proof}
 Manifolds are also locally compact:
 \begin{proposition}\label{prop:nociones-fundamentales-variedad-topologica-localmente-compacta}
 Every topological manifold with or without boundary is locally compact.
 \end{proposition}
\begin{proof}
Let \(p\in M\) and choose a chart
\(\phi\colon U\to\phi(U)\subseteq\mathbb H^n\) around \(p\), where
\(\mathbb H^n=\mathbb R^n\) if \(M\) has no boundary. Since \(\phi(U)\)
is open in \(\mathbb H^n\), there exists \(r>0\) such that
\[
\overline B_{\mathrm{euc}}(\phi(p),r)\cap\mathbb H^n
\subseteq\phi(U).
\]
The preimage
\[
K:=\phi^{-1}\bigl(
\overline B_{\mathrm{euc}}(\phi(p),r)\cap\mathbb H^n\bigr)
\]
is compact, since it is homeomorphic to a closed and bounded subset of
\(\mathbb R^n\), and it contains the relatively open neighborhood
\[
\phi^{-1}\bigl(B_{\mathrm{euc}}(\phi(p),r)\cap\mathbb H^n\bigr)
\]
of \(p\).
\end{proof}
 Another covering property of manifolds is paracompactness, which will be very useful in geometric analysis:
 \begin{theorem}[Manifolds are paracompact]\label{teo: variedad es paracompacta}\index{manifolds are paracompact}
 Let $(M,\tau_{M})$ be a topological manifold with or without boundary. Given an open cover $\mathcal{U}\subseteq\tau_{M}$ of $M$ and any basis $\mathcal{B}$ for $\tau_{M}$, there exists a countable, locally finite open refinement of $\mathcal{U}$ consisting of elements of $\mathcal{B}$.
 \end{theorem}
\begin{proof}
We first construct a compact exhaustion. Since $M$ is locally
compact and second countable, it has a countable basis
$(V_j)_{j\geq1}$ of precompact open sets. Set
$K_0:=\varnothing$ and $K_1:=\overline V_1$. Once $K_j$ has been constructed, each point of the compact set
$K_j\cup\overline V_{j+1}$ has a precompact open neighborhood. A
finite subcover of these neighborhoods has union $W$, and the union of their
closures defines a compact set $K_{j+1}$ satisfying
\[
 K_j\cup\overline V_{j+1}
 \subseteq W\subseteq\operatorname{Int}K_{j+1}.
\]
Recursively, we obtain compact sets $(K_j)_{j\geq1}$ such that
\[
 K_j\subseteq\operatorname{Int}K_{j+1},
 \qquad V_1\cup\cdots\cup V_j\subseteq K_j.
\]
Since the sets $V_j$ cover $M$, we have $\displaystyle M=\bigcup_{j\geq1}K_j$.

Define
\[
 C_0:=K_1,
 \qquad W_0:=\operatorname{Int}K_2,
\]
and, for $j\geq1$,
\[
 C_j:=K_{j+1}\setminus\operatorname{Int}K_j,
 \qquad
 W_j:=\operatorname{Int}K_{j+2}\setminus K_{j-1}.
\]
Each $C_j$ is compact and contained in $W_j$. For $j\geq0$ and
$x\in C_j$, choose $U_x\in\mathcal U$ with $x\in U_x$ and, since
$\mathcal B$ is a basis, an element $B_x\in\mathcal B$ such that
\[
 x\in B_x\subseteq U_x\cap W_j.
\]
Compactness of $C_j$ allows us to retain finitely many of these
open sets. The union of the finite families obtained for
$j=0,1,2,\ldots$ is countable and refines $\mathcal U$.

We show that this family covers $M$. Given $x\in M$, there exists $j\geq0$ such
that $x\in\operatorname{Int}K_{j+1}$; choose the smallest such index.
If $j=0$, then $x\in C_0$. If $j\geq1$, minimality implies
$x\notin\operatorname{Int}K_j$, while $x\in K_{j+1}$, and
therefore $x\in C_j$.

It remains to verify local finiteness. Let $x\in M$ and choose $N\geq1$ with
$x\in\operatorname{Int}K_N$. If $j\geq N+1$, then
$K_N\subseteq K_{j-1}$ and
\[
 \operatorname{Int}K_N\cap W_j=\varnothing.
\]
Thus the neighborhood $\operatorname{Int}K_N$ of $x$ can intersect only
elements arising from the indices $0,\ldots,N$, and for each of
these we chose a finite family. The refinement is locally finite.

The proof uses only the fact that charts provide a basis
of precompact open sets; it therefore applies unchanged to charts taking values in
relatively open subsets of $\mathbb H^n$.
\end{proof}
Various coordinate systems can be defined on a smooth manifold. Among them is the notion of a \textit{regular coordinate ball}, which generalizes the notion of a \textit{ball} in Euclidean spaces to manifolds.
\begin{definition}\label{def:nociones-fundamentales-norma-variedad-suave-label-bola}\index{regular coordinate ball}
 If $M$ is a smooth manifold without boundary, we use the following terminology:
\begin{enumerate}[label=(\alph*)]
 \item
 $B\subseteq M$ is a coordinate ball centered at $x$ (or at $\phi^{-1}(x)$) if there exists a homeomorphism $\phi\colon B\longrightarrow B_{\mathrm{euc}}(x,r)$, where \[B_{\mathrm{euc}}(x,r):=\{y\in\mathbb{R}^{n}\mid \|y-x\|<r\},\] for some $r>0$ and $x\in\mathbb{R}^{n}$, such that $(B,\phi)$ is a smooth chart of $M$.
 \item $B\subseteq M$ is a regular coordinate ball (centered at $x$) if there exists a coordinate ball $B'$ such that $\overline{B}\subseteq B'$ and there is a smooth coordinate map $\phi\colon B'\longrightarrow \mathbb{R}^{n}$ such that, for some positive numbers $r'>r>0$, we have $\phi(B)=B_{\mathrm{euc}}(x,r)$, $\phi(\overline{B})=\overline{B}_{\mathrm{euc}}(x,r)$, and $\phi(B')=B_{\mathrm{euc}}(x,r')$.
\end{enumerate}
If $M$ has boundary, we use the same terminology, replacing Euclidean
balls by the relative balls
$B_{\mathrm{euc}}(x,r)\cap\mathbb H^n$ in a boundary chart and taking
closures in $\mathbb H^n$.
\end{definition}
\begin{figure}[h]
 \centering
 \includegraphics[scale=0.5]{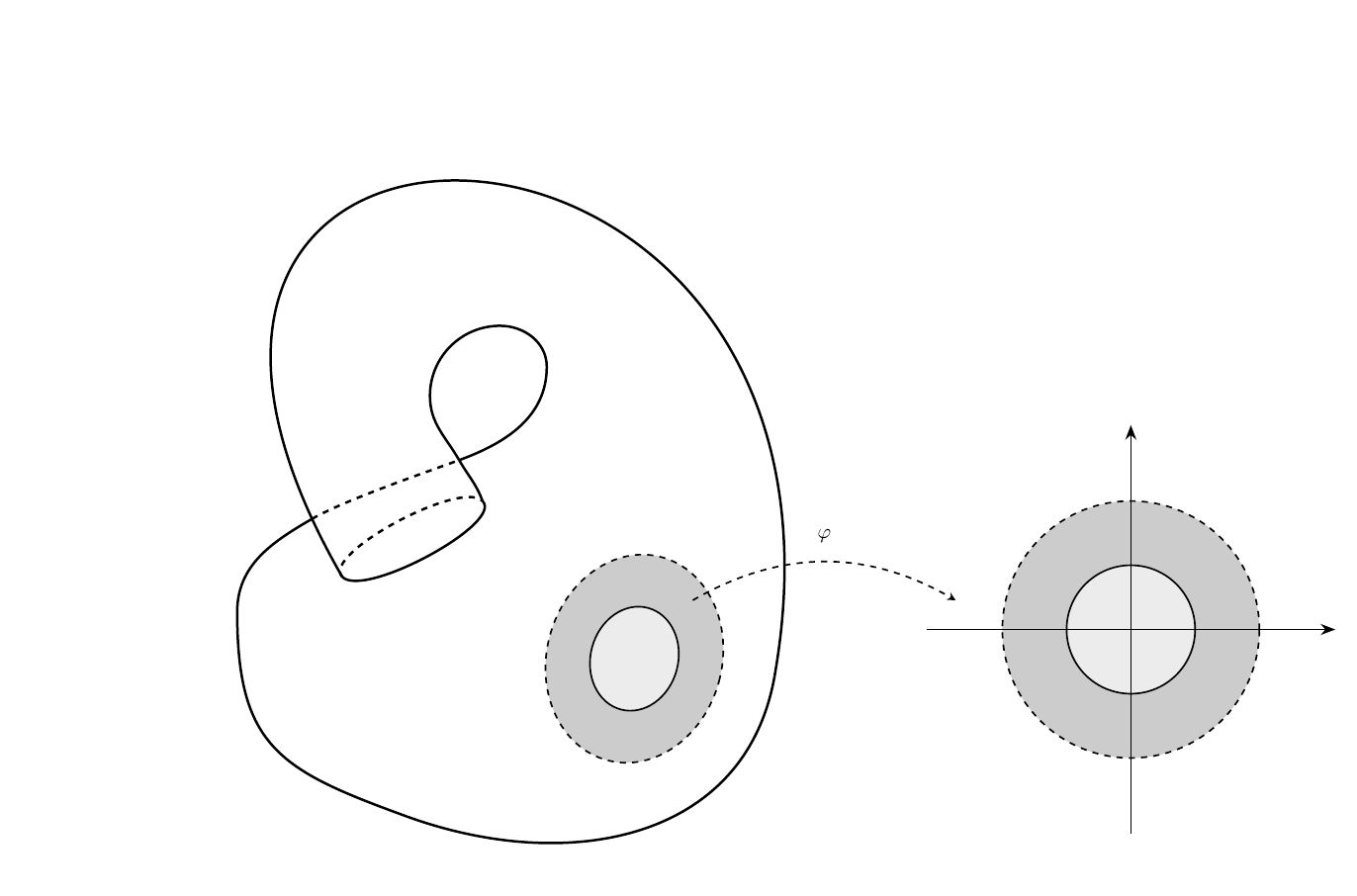}
 \caption{A precompact regular coordinate ball
 whose coordinate map extends to the boundary.}
 \label{fig:bola-coordenada-regular}
\end{figure}

\begin{remark}\label{obs:nociones-fundamentales-bolas-coordenadas-regulares-compactas}
 Every regular coordinate ball $(B,\phi)$ has the property that $\overline{B}$ is compact and $\phi$ extends to $\overline{B}$. We will call domains of this type \textit{\textbf{compact regular coordinate balls}}.
\end{remark}
Regular coordinate balls are important because of their topological properties. Moreover, every smooth manifold admits a countable atlas of regular coordinate balls, as the following proposition shows:
\begin{proposition}
\label{bolacoordenadaregular}
 Every smooth manifold with or without boundary has a countable basis of regular coordinate balls, and therefore admits a countable, locally finite open cover by regular coordinate balls.
\end{proposition}
\begin{proof}
By second countability, we extract a countable subcover from the cover by
chart domains. Within each coordinate image,
we take Euclidean balls, or half-balls in a boundary
chart, with rational centers and radii, whose closures lie in
slightly larger balls in the same coordinate domain. Their preimages
form a countable family of regular coordinate balls and constitute
a basis for $M$. Applying
Theorem~\ref{teo: variedad es paracompacta} to this basis and to the cover
$\{M\}$ gives a countable, locally finite subfamily covering
$M$.
\end{proof}

 Finally, to complete the definitions in this section, we specify the
appropriate model for a manifold with boundary:
\[
 \mathbb H^n:=\mathbb R^{n-1}\times[0,\infty),\qquad
 \operatorname{Int}(\mathbb H^n)=\mathbb R^{n-1}\times(0,\infty),
 \qquad \partial\mathbb H^n=\mathbb R^{n-1}\times\{0\}.
\]
All open subsets of $\mathbb H^n$ are understood in the relative topology.
A function defined on one of them is smooth if, in a neighborhood of each
point, it is the restriction of a smooth function defined on an open subset of
$\mathbb R^n$.

 \begin{definition}\label{def:nociones-fundamentales-variedad-topologica-con-frontera-de-dimension}\index{topological manifold with boundary of dimension@topological manifold with boundary of dimension} Let $M$ be a topological space.

 \begin{enumerate}[label=(\alph*)]
 \item We say that $M$ is a \textbf{topological manifold with boundary of dimension $n$} if it is Hausdorff and second countable, and each point has an open neighborhood homeomorphic to a relatively open subset of $\mathbb H^n$.
 \item A chart for $M$ is a pair $(U,\phi)$, where $U\subseteq M$ is open and $\phi\colon U\longrightarrow \phi(U)\subseteq\mathbb H^n$ is a homeomorphism onto a relatively open subset.
 \item A chart $(U,\phi)$ is an \textbf{interior chart} if $\phi(U)\subseteq\operatorname{Int}(\mathbb H^n)$; in that case, its image is an open subset of $\mathbb R^n$. It is a \textbf{boundary chart} if $\phi(U)\cap\partial\mathbb H^n\neq\varnothing$. The domain of a boundary chart may contain both interior and boundary points.
 \item An atlas for $M$ is a collection of charts $\mathcal A$ such that $\displaystyle\bigcup_{(U,\phi)\in\mathcal A}U=M$.
 \end{enumerate}
 \end{definition}

 Compatibility and maximality for charts of a manifold
with boundary are defined as follows:

 \begin{definition}\label{def:nociones-fundamentales-variedad-3}\index{smooth structure on a manifold with boundary}
 Let $M$ be a topological manifold with boundary of dimension $n$.
 \begin{enumerate}[label=(\alph*)]
 \item We say that two charts $(U,\phi)$ and $(V,\psi)$ are smoothly compatible if $U\cap V= \varnothing$ or if the transition map \[\psi\circ \phi^{-1}\colon \phi(U\cap V)\longrightarrow \psi(U\cap V)\] is a diffeomorphism of class $C^{\infty}$.

 \item A smooth atlas is an atlas $\mathcal{A}$ such that any two charts in $\mathcal{A}$ are smoothly compatible.
 \item A smooth structure on $M$ is a smooth atlas $\mathcal{A}$ maximal with respect to $\subseteq$ among the smooth atlases of $M$.
 \end{enumerate}
 \end{definition}

 \begin{definition}\label{def:nociones-fundamentales-variedad-4}\index{smooth manifold with boundary}
 A smooth manifold with boundary of dimension $n$ is a topological manifold with boundary of dimension $n$ equipped with a smooth structure.
 \end{definition}

\begin{theorem}[Smooth invariance of the boundary]
\label{teo:invariancia-suave-frontera-semiespacio}
\index{invariance of the boundary}
Let $(U,\phi)$ and $(V,\psi)$ be smoothly compatible charts of a
manifold with boundary, and let $p\in U\cap V$. Then
\[
 \phi(p)\in\partial\mathbb H^n
 \quad\Longleftrightarrow\quad
 \psi(p)\in\partial\mathbb H^n.
\]
Consequently, the notions of interior point and boundary point obtained
by examining the image of a point do not depend on the chart.
\end{theorem}

The proof follows the strategy of Lee~\cite{LeeS}.

\begin{proof}
Suppose, by contradiction, that $x:=\phi(p)$ is interior and that
$y:=\psi(p)$ belongs to $\partial\mathbb H^n$. Restricting the charts,
we may take a Euclidean ball $B$ around $x$, contained in
$\phi(U\cap V)\cap\operatorname{Int}(\mathbb H^n)$. The change of
coordinates
\[
 F:=\psi\circ\phi^{-1}\colon B\longrightarrow\mathbb H^n
\]
is smooth, and its inverse admits a smooth extension near $y$ to an
open subset of $\mathbb R^n$. The identity $F^{-1}\circ F=\operatorname{Id}$
on $B$ and the chain rule show that $dF_x$ has a left
inverse; since its domain and codomain have dimension $n$, $dF_x$ is
invertible. Theorem~\ref{lem:funcion-inversa-euclidea-previa-variedades} implies that $F(B)$
contains a Euclidean neighborhood of $y$. This is impossible because
$F(B)\subseteq\mathbb H^n$ and no Euclidean neighborhood of a point of
$\partial\mathbb H^n$ is contained in $\mathbb H^n$. Interchanging the charts rules out the opposite possibility.
\end{proof}

\begin{definition}\label{def:puntos-interiores-frontera-cartas}
\index{interior point of a manifold with boundary}
\index{boundary point of a manifold}
Let $M$ be a smooth manifold with boundary. A point $p\in M$ is
an \textbf{interior point} if $\phi(p)\in\operatorname{Int}(\mathbb H^n)$ in one
---and, by the preceding theorem, in every--- chart containing it. It is a
\textbf{boundary point} if
$\phi(p)\in\partial\mathbb H^n$. We denote the corresponding sets
by $\operatorname{Int}(M)$ and $\partial M$.

A boundary chart is \textbf{centered at $p\in\partial M$} if,
after a translation in the first $n-1$ coordinates,
$\phi(p)=0$. This condition concerns the image of the distinguished point;
it does not mean that the entire chart domain consists of boundary
points.
\end{definition}

 \begin{remark}[Convention for finite-dimensional manifolds]
 Throughout the finite-dimensional treatment, the term manifold includes
 the Hausdorff and second countability assumptions in the preceding definitions;
 we will not repeat them in every statement.
 \end{remark}

 \begin{remark}\label{obs:nociones-fundamentales-conveniente-hacer-notar-siguiente-frontera-variedad}
 The following properties hold:
 \begin{enumerate}
 \item
 The boundary $\partial M$ and the interior $\operatorname{Int}(M)$ are disjoint and their union is $M$. Moreover, $\operatorname{Int}(M)$ is open and $\partial M$ is closed.
 \item
 A smooth manifold may be regarded as the special case of a smooth manifold with empty boundary, or one in which every point of $M$ is an interior point.
 \end{enumerate}

 \end{remark}

\subsection{Examples of manifolds with boundary}

\begin{example}[A closed interval]
\label{ej:variedades-intervalo-cerrado-frontera}
Let $a<b$ and $\ell:=b-a$. The space $I:=[a,b]$ is compact, Hausdorff,
and second countable. The relatively open sets
\[
 U_a:=[a,a+2\ell/3),\qquad U_b:=(a+\ell/3,b]
\]
cover it. The charts
\[
 \phi_a(x):=x-a,\qquad \phi_b(x):=b-x
\]
have image $[0,2\ell/3)$, open in $\mathbb H^1=[0,\infty)$.
They are homeomorphisms, and the transition on the overlap is
$t\mapsto\ell-t$, with the same expression for its inverse. This map
extends smoothly to $\mathbb R$. Hence $I$ is a smooth manifold of
dimension one with
\[
 \operatorname{Int}(I)=(a,b),\qquad \partial I=\{a,b\}.
\]
At the right endpoint we use $b-x$ so that the normal coordinate is
nonnegative. A smooth function on $I$ may be nonzero at its endpoints;
for example, the constant function one is smooth.
\end{example}

\begin{example}[The closed ball]
\label{ej:variedades-bola-cerrada-frontera}
For $n\geq2$, let
$\overline B^n:=\{x\in\mathbb R^n\mid |x|\leq1\}$, with its subspace
topology. It is compact, Hausdorff, and second countable. Its open ball
is an interior chart domain with the identity. We construct the
remaining charts around each $p\in S^{n-1}$.

An orthogonal transformation sends $p$ to $(0,\ldots,0,1)$. Work in
these coordinates and fix $0<c<1$. On the relatively open cap
\[
 U_c:=\{x\in\overline B^n\mid x^n>c\}
\]
define, writing $x=(x',x^n)$,
\[
 \phi_c(x):=\bigl(x',\sqrt{1-|x'|^2}-x^n\bigr)=(y,t).
\]
Since $x^n>0$ and $|x|\leq1$, we have
$x^n\leq\sqrt{1-|x'|^2}$; hence $t\geq0$. The image is exactly
\[
 \Omega_c:=\left\{(y,t)\in\mathbb R^{n-1}\times[0,\infty)
 \middle| |y|^2<1-c^2,\quad
                t<\sqrt{1-|y|^2}-c\right\}.
\]
It is relatively open in $\mathbb H^n$. The inverse is
\[
 \phi_c^{-1}(y,t)=\bigl(y,\sqrt{1-|y|^2}-t\bigr).
\]
The inequalities defining $\Omega_c$ ensure that the last component
is greater than $c$ and at most $\sqrt{1-|y|^2}$, so the inverse takes
values in $U_c$. Both expressions are continuous and are restrictions
of smooth maps on Euclidean open sets where $|y|<1$. If caps are chosen
using different orthogonal transformations, the transitions are
compositions of these expressions with orthogonal linear maps; they
and their inverses admit local smooth extensions. These caps, together
with the interior chart, cover the closed ball and form a smooth
boundary atlas.

In each cap, $t=0$ is equivalent to $|x|=1$, whereas $t>0$ is
equivalent to $|x|<1$. Thus
\[
 \operatorname{Int}(\overline B^n)=B^n,\qquad
 \partial\overline B^n=S^{n-1}.
\]
This description obtains the boundary structure directly from the
charts. For $n=1$ the closed ball is the interval $[-1,1]$ of the
preceding example. A cap chart is not the projection onto $x'$:
that projection loses the information about the normal distance to
the face; coordinate $t$ retains it.
\end{example}

 \begin{figure}[h]
 \includegraphics[scale=0.2]{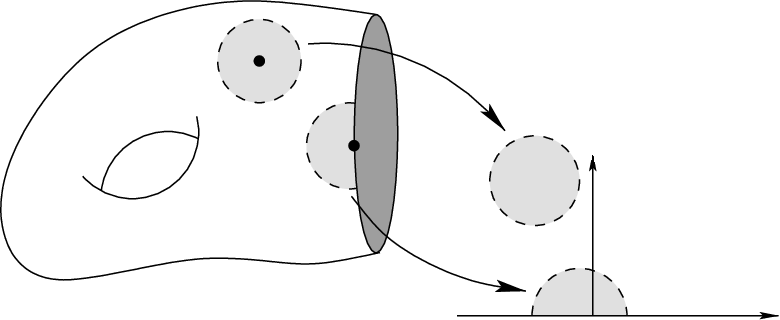}
 \centering
 \label{fig:frontera-variedad-suave}
 \end{figure}

The abstract definition is especially useful when a space is constructed by gluing local models. The following lemma collects the conditions under which a family of coordinates uniquely determines the topology and smooth structure of the resulting set.
\begin{lemma}[Smooth chart lemma]\label{lema de la carta suave}\index{lemma!smooth chart}
 Let $M$ be a set, and suppose we have a collection $\{U_{\alpha}\}_{\alpha\in A}$ of subsets of $M$ and maps $\phi_{\alpha}\colon U_{\alpha}\longrightarrow \mathbb{R}^{n}$ (or $\phi_{\alpha}\colon U_{\alpha}\longrightarrow \mathbb{H}^{n}$) such that:
 \begin{enumerate}[label=(\alph*)]
 \item For each $\alpha \in A$, $\phi_{\alpha}$ is a bijection between $U_{\alpha}$ and an open subset $\phi_{\alpha}(U_{\alpha})\subseteq\mathbb{R}^{n}$ (or a relatively open subset $\phi_{\alpha}(U_{\alpha})\subseteq\mathbb{H}^{n}$).
 \item For each $\alpha,\beta\in A$, the sets $\phi_{\alpha}(U_{\alpha}\cap U_{\beta})$ and $\phi_{\beta}(U_{\alpha}\cap U_{\beta})$ are open in $\mathbb{R}^{n}$ (or relatively open in $\mathbb{H}^{n}$).
 \item If $U_{\alpha}\cap U_{\beta}\neq \varnothing$, the map $\phi_{\beta}\circ \phi_{\alpha}^{-1}\colon \phi_{\alpha}(U_{\alpha}\cap U_{\beta})\longrightarrow \phi_{\beta}(U_{\alpha}\cap U_{\beta})$ is smooth.
 \item There exists a countable $N\subseteq A$ such that $M=\displaystyle\bigcup_{\alpha\in N}U_{\alpha}$.
 \item If $p,q\in M$ and $p\neq q$, then either there exists $\alpha\in A$ such that $p,q\in U_{\alpha}$, or there exist $\alpha,\beta\in A$ such that $U_{\alpha}\cap U_{\beta}=\varnothing$ with $p\in U_{\alpha}$ and $q\in U_{\beta}$.
 \end{enumerate}
 Then $M$ has a unique topology and a unique smooth structure such that $M$ is a smooth manifold (with boundary) and $(U_{\alpha},\phi_{\alpha})$ is a smooth chart for every $\alpha\in A$.
\end{lemma}
\begin{proof}
We write $\mathbb E^n$ for the appropriate model space,
$\mathbb R^n$ or $\mathbb H^n$, always using the relative topology
in the latter case. Consider the family
\[
 \mathcal C:=
 \left\{\phi_\alpha^{-1}(V)\middle|
 \alpha\in A,\ V\subseteq\phi_\alpha(U_\alpha)
 \text{ open}\right\}.
\]
This family covers $M$. If $p$ belongs to
$\phi_\alpha^{-1}(V)\cap\phi_\beta^{-1}(W)$, then
\[
 \phi_\alpha^{-1}(V)\cap\phi_\beta^{-1}(W)
 =\phi_\alpha^{-1}\!\left(
 V\cap(\phi_\beta\circ\phi_\alpha^{-1})^{-1}(W)
 \right).
\]
The set in parentheses is open in
$\phi_\alpha(U_\alpha)$: we first restrict to the open set
$\phi_\alpha(U_\alpha\cap U_\beta)$ and then use continuity of the
change of coordinates. Therefore, $\mathcal C$ is a basis for a
topology $\tau$ on $M$.

By construction, each $U_\alpha$ is open and $\phi_\alpha$ is a
homeomorphism from $U_\alpha$ onto its image. We show that $(M,\tau)$ is
Hausdorff. If two distinct points $p,q$ belong to the same
$U_\alpha$, we separate their images in the Hausdorff space
$\mathbb E^n$ and take preimages. If we do not use a common chart, the
last assumption provides disjoint domains $U_\alpha$ and $U_\beta$
containing $p$ and $q$, respectively. In both cases we obtain
disjoint open neighborhoods.

Let $N\subseteq A$ be countable, as in the fourth assumption. Each
$\phi_\alpha(U_\alpha)$, with $\alpha\in N$, has a countable basis; the
preimages of all these bases form a countable family. This is a basis
for $M$, since the sets $U_\alpha$, $\alpha\in N$, cover the space. Thus
$(M,\tau)$ is second countable and locally Euclidean, or
locally modeled on $\mathbb H^n$.

For $U_\alpha\cap U_\beta\neq\varnothing$, the third assumption, together with the
same assumption with $\alpha$ and $\beta$ interchanged, shows that
$\phi_\beta\circ\phi_\alpha^{-1}$ is a diffeomorphism between the corresponding
open sets. The given charts therefore form a smooth atlas. The maximal
atlas containing it defines a smooth structure on $M$. In the case of
$\mathbb H^n$, the same construction yields boundary charts, and all
local extensions involved in the notion of smoothness are taken
on open subsets of $\mathbb R^n$; thus the argument also covers
manifolds with boundary.

Finally, any other topology making all the maps $\phi_\alpha$
homeomorphisms onto their images would necessarily have $\mathcal C$
as a basis, and would therefore coincide with $\tau$. A smooth structure for
which the given charts were smooth would have the same maximal atlas,
since a chart belongs to this atlas if and only if it is smoothly compatible
with every chart $\phi_\alpha$. This proves both uniqueness assertions.
\end{proof}

\subsection{Manifolds with corners, depth, and products}

The product of two closed half-spaces has two nonnegative normal
coordinates. This model requires an enlargement of boundary charts.
The enlargement retains smoothness through local extensions to
Euclidean open sets.

\begin{definition}[Models and manifolds with corners]
\label{def:variedades-con-esquinas-modelos}
For $n\in\mathbb N_0$ and $r\in\{0,\ldots,n\}$, define
\[
 G_{n,r}:=\mathbb R^{n-r}\times[0,\infty)^r.
\]
A map between relatively open subsets of these models is of class
$C^k$, with $k\in\mathbb N_0\cup\{\infty\}$, if near every point its
components are restrictions of functions of class $C^k$ on a Euclidean
open set. It is a diffeomorphism if it is bijective and both it and its
inverse are smooth. A smooth manifold with corners of dimension $n$
is a Hausdorff, second-countable space equipped with a maximal atlas of
homeomorphisms $\phi\colon U\to\Omega$, where
$\Omega\subseteq G_{n,r}$ is relatively open and $r$ may depend on the
chart, whose transitions are diffeomorphisms in this sense. A map
between manifolds with corners is of class $C^k$ if its chart
representations are of class $C^k$.
\end{definition}

Composition is well defined: an extension of the first representation
takes, on a sufficiently small neighborhood of the point, values in
the open domain of an extension of the second. Their composition is
an extension of class $C^k$ of the composite map. This also proves
independence of charts in the definition of a map. For $r=0$ we recover
the Euclidean model, and for $r=1$ the half-space. Smooth maps are not
required to preserve faces; for instance, $t\mapsto t^2$ on
$[0,\infty)$ is smooth. Additional conditions on maps between faces
are not part of this definition. For $k=0$, the condition is ordinary
continuity: a continuous representation extends locally by replacing
each normal coordinate by its absolute value and restricting to a
small ball whose image under this replacement stays in the relative
domain.

\begin{lemma}[Uniqueness of jets and smoothness up to the faces]
\label{lem:esquinas-jets-extension-local}
Let $\Omega\subseteq G_{n,r}$ be relatively open. Two extensions of
class $C^k$ of the same function on $\Omega$ have the same derivatives
of order at most $k$ at all points of $\Omega$. Moreover, a function
is smooth in the sense of
Definition~\ref{def:variedades-con-esquinas-modelos} if and only if it
is smooth on $\Omega\cap\operatorname{Int}_{\mathbb R^n}G_{n,r}$ and
all its ordinary derivatives extend continuously to $\Omega$.
\end{lemma}

\begin{proof}
Every point of $\Omega$ is a limit of points in its Euclidean interior:
increase the zero normal coordinates slightly while keeping the
others fixed, inside a relative ball contained in $\Omega$.
Two extensions agree on that open interior; their derivatives agree
there. Continuity of derivatives up to order $k$ gives equality at
points of the faces. This includes order zero.

A smooth extension, by restriction, has the asserted continuity of
all derivatives. Conversely, fix $z\in\Omega$ and choose a Euclidean
cutoff $\chi$ supported in a ball whose intersection with $G_{n,r}$
is contained in $\Omega$, and equal to one near $z$.
Extend $\chi f$ by zero to $G_{n,r}$; it is smooth in the interior and
its derivatives extend continuously to the faces, because its support
does not reach the edge of the chosen relative domain.
Lemma~\ref{lem:extension-suave-semiespacio-soporte-fijo} of Appendix~B
extends this function to $\mathbb R^n$. The extension agrees with $f$
near $z$. Since $z$ was arbitrary, we obtain the required local
extensions.
\end{proof}

\begin{proposition}[Invariance of depth]
\label{prop:esquinas-profundidad-invariante}
In a chart $\phi(p)=z\in G_{n,r}$, the number
\[
 d_M(p):=\#\{i\in\{n-r+1,\ldots,n\}\mid z^i=0\}
\]
is independent of the chart. It is called the depth of $p$. The sets
\[
 S_j(M):=\{p\in M\mid d_M(p)=j\},\qquad 0\leq j\leq n,
\]
are smooth manifolds without boundary of dimension $n-j$ when nonempty.
The interior $\operatorname{Int}(M):=S_0(M)$ is open and dense, and
\[
 \partial M:=M\setminus \operatorname{Int}(M)
           =\bigcup_{j=1}^n S_j(M)
\]
is closed. Here $\partial M$ denotes the subset of points of positive
depth, with its subspace topology.
\end{proposition}

\begin{proof}
Let $F\colon\Omega\to\Omega'$ be a coordinate transition and let
$H:=F^{-1}$. Take smooth extensions $\widetilde F$ and
$\widetilde H$ near $z$ and $F(z)$. After shrinking their domains,
$\widetilde H\circ\widetilde F$ is defined and agrees with the
identity on a relative neighborhood of $z$. Uniqueness of jets in
Lemma~\ref{lem:esquinas-jets-extension-local} implies
\[
 D\widetilde H(F(z))D\widetilde F(z)=I_n.
\]
The same argument for the other composition gives the identity in the
opposite order. Thus $D\widetilde F(z)$ is invertible even on the
faces. We henceforth write it as $DF(z)$.

Consider the cone
\[
 C_z:=\{v\in\mathbb R^n\mid
        v^i\geq0\ \text{if }i>n-r\text{ and }z^i=0\}.
\]
It is exactly the set of velocities at zero of curves of class $C^1$
defined on $[0,\varepsilon)$, starting at $z$ and taking values in
$\Omega$. Indeed, if the initial coordinate is zero and the curve
remains nonnegative, its right derivative is nonnegative. Conversely,
$z+tv$ remains in $G_{n,r}$ for sufficiently small $t\geq0$:
initially positive normal coordinates stay positive, and zero ones
grow with nonnegative velocity. Relative openness ensures that it
also remains in $\Omega$.
Applying the chain rule to these curves and $F$ gives
$DF(z)C_z\subseteq C_{F(z)}$; applying it to $H$ gives the reverse
inclusion. Consequently,
\[
 DF(z)C_z=C_{F(z)},\qquad
 DF(z)(C_z\cap(-C_z))=C_{F(z)}\cap(-C_{F(z)}).
\]
The subspace $C_z\cap(-C_z)$ is defined by setting precisely the
$d_M(p)$ active normal coordinates to zero, and has dimension
$n-d_M(p)$. Since a linear isomorphism preserves this dimension,
depth is invariant.

Around a point of depth $j$, shrink the chart so that initially
positive normal coordinates remain positive. They can then be treated
as free coordinates. A permutation and a translation of the free
coordinates give a chart modeled on $G_{n,j}$, whose $j$ normal
coordinates vanish at the point. In this neighborhood, points of
depth $j$ are exactly those where these $j$ coordinates vanish.
Restricting charts and deleting the zero coordinates gives images
open in $\mathbb R^{n-j}$. The restricted transitions and their
inverses are smooth, since they arise from smooth extensions of the
original transitions. The subspace topology remains Hausdorff and
second countable. This proves the assertion about $S_j(M)$. Points
where all normal coordinates are positive form an open dense subset
of every chart, proving the remaining assertions.
\end{proof}

\begin{remark}[The total boundary and faces]
\label{obs:esquinas-frontera-total-caras}
At a point of depth one the total boundary is described by one
coordinate equation and is locally smooth. At a point of depth
$j\geq2$ it is the union of $j$ intersecting coordinate faces.
This definition does not assign it the structure of a smooth manifold
without boundary of dimension $n-1$. Each coordinate face is treated
separately; intersections give strata of greater depth. Nor do we
assume that the global closure of each component of $S_1(M)$ is an
embedded face: that is an additional hypothesis on some classes of
manifolds with corners. Our notation $\partial M$ retains the boundary
as a subset; constructions separating the distinct faces over a single
point require another definition.
\end{remark}

\begin{theorem}[Product of manifolds with corners]
\label{teo:producto-variedades-con-esquinas}
Let $M$ and $N$ be smooth manifolds with corners of dimensions $m$ and
$n$. The product topology and product charts determine a canonical
smooth structure with corners on $M\times N$, of dimension $m+n$.
The projections are smooth. A map from a manifold with corners $P$
to $M\times N$ is of class $C^k$ if and only if its two components
are of class $C^k$. Moreover,
\[
 \begin{aligned}
 d_{M\times N}(p,q)&=d_M(p)+d_N(q),\\
 S_j(M\times N)&=
   \coprod_{a+b=j} S_a(M)\times S_b(N),\\
 \operatorname{Int}(M\times N)&=\operatorname{Int}(M)\times \operatorname{Int}(N),\\
 \partial(M\times N)&=(\partial M\times N)\cup(M\times\partial N).
 \end{aligned}
\]
The union of strata is disjoint; the two parts in the last union
intersect in $\partial M\times\partial N$.
\end{theorem}

\begin{proof}
The product is Hausdorff: two distinct pairs differ in at least one
factor, and disjoint open sets in that factor give disjoint product
open sets. If $(B_i)$ and $(C_j)$ are countable bases of the factors,
$(B_i\times C_j)_{i,j\in\mathbb N}$ is a countable basis of the product.

Let $\phi\colon U\to\Omega\subseteq G_{m,r}$ and
$\psi\colon V\to\Lambda\subseteq G_{n,s}$ be charts. Write their
coordinates as $(x',t)$ and $(y',u)$, where $t$ has $r$ components
and $u$ has $s$. The linear permutation
\[
 P_{r,s}\bigl((x',t),(y',u)\bigr):=(x',y',t,u)
\]
identifies $G_{m,r}\times G_{n,s}$ with $G_{m+n,r+s}$ and is a
homeomorphism. Hence
\[
 \Phi:=P_{r,s}\circ(\phi\times\psi)\colon
 U\times V\longrightarrow P_{r,s}(\Omega\times\Lambda)
\]
is a homeomorphism onto a relatively open subset of this model.
These domains cover $M\times N$.

If $\widehat\phi$ and $\widehat\psi$ are other charts and
$\widehat P$ is the corresponding permutation, the product transition is
\[
 \widehat P\circ
 \bigl((\widehat\phi\circ\phi^{-1})\times
       (\widehat\psi\circ\psi^{-1})\bigr)\circ P_{r,s}^{-1}.
\]
Near every point, smooth extensions of the factor transitions give,
by taking their product and the two linear permutations, a smooth
extension of this transition. Its inverse extends in the same way.
Thus we have a compatible atlas. Adding all compatible charts gives
its unique maximal atlas. Replacing factor atlases by compatible
atlases gives transitions of the same form, so does not change the
product structure.

In product coordinates, the projections are linear projections after
the inverse permutation; they are smooth. If the components of a map
from $P$ are of class $C^k$, their local extensions can be chosen on
one common smaller neighborhood, and combining them gives an extension
of the product map. The converse follows by composition with the
smooth projections. This proves the criterion for maps.

Finally, the zero normal coordinates in the product chart are those
of $\phi(p)$ together with those of $\psi(q)$. Their number is the
sum of the depths. The established invariance turns this coordinate
count into an intrinsic identity. The stratum, interior, and boundary
formulas follow. A pair lies in both parts of the boundary exactly
when both components are boundary points, which gives the asserted
intersection.
\end{proof}

\begin{corollary}[Two factors with smooth boundary]
\label{cor:producto-dos-variedades-frontera}
If $M$ and $N$ are smooth manifolds with boundary, their product has
depth at most two. Points in $\partial M\times\partial N$ have depth
two. If this set is nonempty, the product structure admits no
compatible atlas consisting exclusively of Euclidean and half-space
charts. If at most one factor has nonempty boundary, the product is
a smooth manifold with boundary in the previous sense.
\end{corollary}

\begin{proof}
In each factor, depth is zero or one. The theorem gives their sum.
A half-space chart has depth at most one, so invariance of depth
excludes a compatible chart of this type at a point of depth two.
If there are no points of depth two, the local reduction to
$G_{m+n,j}$ in the proof of
Proposition~\ref{prop:esquinas-profundidad-invariante} uses only $j=0,1$,
and provides the required boundary atlas.
\end{proof}

\begin{example}[Tori, cylinders, and annuli]
\label{ej:variedades-toros-cilindros-anillos}
The torus $\mathbb T^n:=(S^1)^n$ is a compact manifold without boundary
of dimension $n$: repeatedly apply the product theorem to the charts
of $S^1$ in Example~\ref{ej:variedades-esfera-cartas-hemisferios}.
Compactness follows from the finite product of compact spaces, and
depth is zero.

For $n\geq2$ and $0<a<b$, the cylinder $S^{n-1}\times[a,b]$ is a
compact manifold with boundary
\[
 \partial(S^{n-1}\times[a,b])
     =(S^{n-1}\times\{a\})\coprod(S^{n-1}\times\{b\}).
\]
Its charts combine a sphere chart with the interval coordinate $r-a$
or $b-r$. The map
\[
 F\colon S^{n-1}\times[a,b]\longrightarrow
 A_{a,b}:=\{x\in\mathbb R^n\mid a\leq|x|\leq b\},
 \qquad F(\omega,r)=r\omega
\]
is a homeomorphism, with inverse $x\mapsto(x/|x|,|x|)$.
The inverse is continuous because $|x|\geq a>0$. Transporting the
product charts by $F$ gives a boundary atlas on the annulus $A_{a,b}$.
Its expressions and those of $F$ are smooth, since the norm is smooth
off the origin and the sphere charts have smooth inverses. The boundary
consists of the two spheres of radii $a$ and $b$. This constructs,
rather than merely asserts, the smooth boundary structure of the
annulus.
\end{example}

\begin{example}[Charts and strata of the square and the cube]
\label{ej:variedades-cuadrado-cubo-esquinas}
\index{manifold with corners!square and cube}
Let $n\geq1$ and let $Q^n=[0,1]^n$, with the topology induced
from $\mathbb R^n$. It is Hausdorff and second countable as a
Euclidean subspace, and compact because it is closed and bounded.
We construct the charts of its product structure explicitly.

Fix $p\in Q^n$ and set
\[
 J_0:=\{i\mid p^i=0\},\qquad J_1:=\{i\mid p^i=1\},\qquad
 F:=\{i\mid0<p^i<1\},\qquad d:=|J_0|+|J_1|.
\]
Choose $0<\varepsilon<1/2$ smaller than $p^i$ and $1-p^i$ for
every $i\in F$; if $F$ is empty, only the first condition is needed.
Let $U_p=\prod_{i=1}^n I_i$, where
\[
 I_i:=\begin{cases}
 [0,\varepsilon),&i\in J_0,\\
 (1-\varepsilon,1],&i\in J_1,\\
 (p^i-\varepsilon,p^i+\varepsilon),&i\in F.
 \end{cases}
\]
This set is relatively open and contains $p$. Ordering the indices
of each group, define
\[
 \phi_p(x):=\bigl((x^i-p^i)_{i\in F},
                      (x^i)_{i\in J_0},(1-x^i)_{i\in J_1}\bigr).
\]
It is a homeomorphism from $U_p$ onto
$(-\varepsilon,\varepsilon)^{n-d}\times[0,\varepsilon)^d$,
a relatively open subset of $G_{n,d}$, and $\phi_p(p)=0$.
Its inverse adds $p^i$ in the free coordinates, preserves those
of $J_0$, and subtracts those of $J_1$ from one, restoring the
original order. Both expressions are restrictions of invertible
affine maps on $\mathbb R^n$. Hence every transition
$\phi_q\circ\phi_p^{-1}$ and its inverse are restrictions of
invertible affine maps, even when the two charts use different
depths. Their domains cover $Q^n$ and determine a smooth atlas
with corners compatible with the product charts of the interval.

These charts give
\[
 d_{Q^n}(x)=\#\{i\in\{1,\ldots,n\}\mid x^i\in\{0,1\}\},
 \qquad \operatorname{Int}(Q^n)=(0,1)^n.
\]
For each $j\in\{0,\ldots,n\}$, the stratum $S_j(Q^n)$ is the
disjoint union of $2^j\binom{n}{j}$ open faces of dimension $n-j$:
choose the $j$ fixed coordinates and, for each of them, zero or
one as its value; the others range over $(0,1)$.
Each piece has the preceding free charts and is diffeomorphic
to $(0,1)^{n-j}$. For $j=n$ this model means a point.

Thus $Q^2$ has an interior stratum, four sides without vertices,
and four vertices. In $Q^3$ there are an interior, six open faces,
twelve edges without endpoints, and eight vertices, of depths
$0,1,2,3$, respectively. For example, near $(0,1,0)$ the coordinates
are $(x^1,1-x^2,x^3)$ and the model is $G_{3,3}$.
Invariance of depth prevents replacing a vertex chart by a
compatible half-space chart. This assertion concerns the smooth
product structure: a homeomorphism with a closed ball can
transport another structure, but that does not make it a
diffeomorphism for this atlas.
\end{example}

\begin{example}[The closed simplex and its barycentric coordinates]
\label{ej:variedad-esquinas-simplex-cartas}
\index{manifold with corners!simplex}
For $n\geq1$ consider
\[
 \Delta^n:=\left\{x\in\mathbb R^n\ \middle|\
       x^i\geq0\ (1\leq i\leq n),\quad\displaystyle\sum_{i=1}^n x^i\leq1\right\}.
\]
It is compact, Hausdorff, and second countable. Define the affine
functions $\lambda_0(x)=1-\displaystyle\sum_{i=1}^n x^i$ and $\lambda_i(x)=x^i$
for $1\leq i\leq n$. On $\Delta^n$ they are nonnegative and sum
to one. For $p\in\Delta^n$ let
\[
 A:=\{i\in\{0,\ldots,n\}\mid\lambda_i(p)=0\},\qquad d:=|A|.
\]
At least one coordinate is positive, so $d\leq n$.
Choose $b\notin A$ and write
$F:=\{0,\ldots,n\}\setminus(A\cup\{b\})$.
The domain
\[
 U_p:=\{x\in\Delta^n\mid\lambda_i(x)>0\text{ for every }i\notin A\}
\]
is relatively open and contains $p$. The chart is
\[
 \phi_p(x):=\bigl((\lambda_i(x)-\lambda_i(p))_{i\in F},
                                      (\lambda_i(x))_{i\in A}\bigr).
\]
Writing its coordinates as $(y,t)$, its image is exactly the
following relatively open subset of $G_{n,d}$:
\[
 \Omega_p:=\left\{(y,t)\in\mathbb R^{n-d}\times[0,\infty)^d\ \middle|\
 \begin{array}{l}
 \lambda_i(p)+y_i>0\quad(i\in F),\\
 \lambda_b(p)-\displaystyle\sum_{i\in F}y_i-\displaystyle\sum_{i\in A}t_i>0
 \end{array}\right\}.
\]
Indeed, the inverse reconstructs
\[
 \lambda_i=\lambda_i(p)+y_i\ (i\in F),\qquad
 \lambda_i=t_i\ (i\in A),\qquad
 \lambda_b=\lambda_b(p)-\displaystyle\sum_{i\in F}y_i-\displaystyle\sum_{i\in A}t_i,
\]
and then takes $x^i=\lambda_i$ for $1\leq i\leq n$.
These coordinates are nonnegative, sum to one, and give a point
of $U_p$. Conversely, every point of $U_p$ gives those inequalities.
The chart and its inverse are restrictions of invertible affine
maps between Euclidean spaces of dimension $n$, because the $n$
barycentric coordinates other than $\lambda_b$ determine the
remaining one. They are therefore homeomorphisms, and every
transition and its inverse extend to smooth affine maps.
This constructs the atlas with corners.

In particular,
\[
 d_{\Delta^n}(x)=\#\{i\in\{0,\ldots,n\}\mid\lambda_i(x)=0\},
 \qquad
 \operatorname{Int}(\Delta^n)=\{x\mid\lambda_i(x)>0\ (0\leq i\leq n)\}.
\]
For $0\leq j\leq n$, each choice of $j$ zero coordinates
describes a piece of $S_j(\Delta^n)$ when all the remaining
coordinates are positive. There are $\binom{n+1}{j}$ pieces;
the preceding free coordinates identify each with the interior
of a simplex of dimension $n-j$. The $n+1$ vertices have depth $n$.
For $n=2$ the closed triangle has three sides and three vertices;
for $n=3$ the tetrahedron has four faces, six edges, and four vertices.
Thus a rectangular domain is not required for the corner model:
barycentric charts transform the meeting faces into coordinate faces.
\end{example}

\begin{example}[A region with curved faces and corners]
\label{ej:variedad-esquinas-sector-anular}
\index{manifold with corners!annular sector}
Let $0<a<b$ and
\[
 A:=\{(x,y)\in\mathbb R^2\mid
           a\leq\sqrt{x^2+y^2}\leq b,\quad x\geq0,\quad y\geq0\}.
\]
Define
\[
 F\colon[a,b]\times[0,\pi/2]\longrightarrow A,
 \qquad F(r,\theta)=(r\cos\theta,r\sin\theta).
\]
Every point of $A$ has positive radius and a unique angle in
$[0,\pi/2]$, so $F$ is bijective. It is continuous, its domain
is compact, and $A$ is Hausdorff; hence it is a homeomorphism.
Theorem~\ref{teo:producto-variedades-con-esquinas} gives the
product structure of the rectangle. Transporting its charts
by $F$ gives an atlas with corners on $A$.

We also verify regularity of the expressions, relating this
atlas to the curves visible in the plane. The formula for $F$
extends smoothly to $\mathbb R^2$, and
\[
 DF(r,\theta)=
 \begin{pmatrix}\cos\theta&-r\sin\theta\\
                 \sin\theta&r\cos\theta\end{pmatrix},
 \qquad \det DF(r,\theta)=r>0.
\]
The inverse has radial component $r=\sqrt{x^2+y^2}$.
Where $x>0$ use $\theta=\arctan(y/x)$; where $y>0$ use
$\theta=\pi/2-\arctan(x/y)$. These are smooth expressions on
Euclidean open sets and agree when $x,y>0$, because they describe
the unique first-quadrant angle with that ratio. These two open
sets cover $A$, since the origin does not belong to the region.
Thus the transported charts and their inverses have the required
smooth extensions, including at the endpoints of the arcs.

The depth is
\[
 d_A(F(r,\theta))
 =\mathbf1_{\{a,b\}}(r)+\mathbf1_{\{0,\pi/2\}}(\theta).
\]
The interior is given by $a<r<b$ and $0<\theta<\pi/2$.
The depth-one stratum has four components: two circular arcs
without endpoints and two radial segments without endpoints.
The depth-two points are $(a,0)$, $(b,0)$, $(0,a)$, and $(0,b)$.
For example, near $(a,0)$ the normal coordinates are $(r-a,\theta)$,
and near $(0,b)$ they are $(b-r,\pi/2-\theta)$.
In both cases the local image is a relatively open subset of
$[0,\infty)^2$. Faces can be curved in the ambient space;
what matters is that their intersections are described by
several nonnegative coordinates in a smooth chart.
\end{example}

\section{Smooth maps}

A smooth structure allows us not only to discuss functions on a manifold, but also to compare manifolds through maps that respect calculus. The condition is checked in coordinates, but its validity must be independent of the chosen charts.

Given charts around a point and its image, the map temporarily becomes a function between Euclidean open sets. We require this coordinate representation to be smooth. Compatibility of the atlases ensures that the definition does not depend on the pair of charts used.
 \begin{definition}\label{definicion funcion suave - C^k}\index{smooth map}
 Let $M$ and $N$ be smooth manifolds with or without boundary, and let $F\colon M\longrightarrow N$ be a function. We say that $F$ is a smooth map if, for every $p\in M$, there exist coordinate charts $(U,\phi)$ and $(V,\psi)$ containing $p$ and $F(p)$, respectively, such that $F(U)\subseteq V$ and the map $\psi\circ F\circ \phi^{-1}$ is a smooth function from $\phi(U)$ to $\psi(V)$.

 We say that the map is of class $C^{k}$ if $\psi\circ F\circ \phi^{-1}$ is of class $C^{k}$.
 \end{definition}
 \begin{figure}[h]
\includegraphics[scale=0.4]{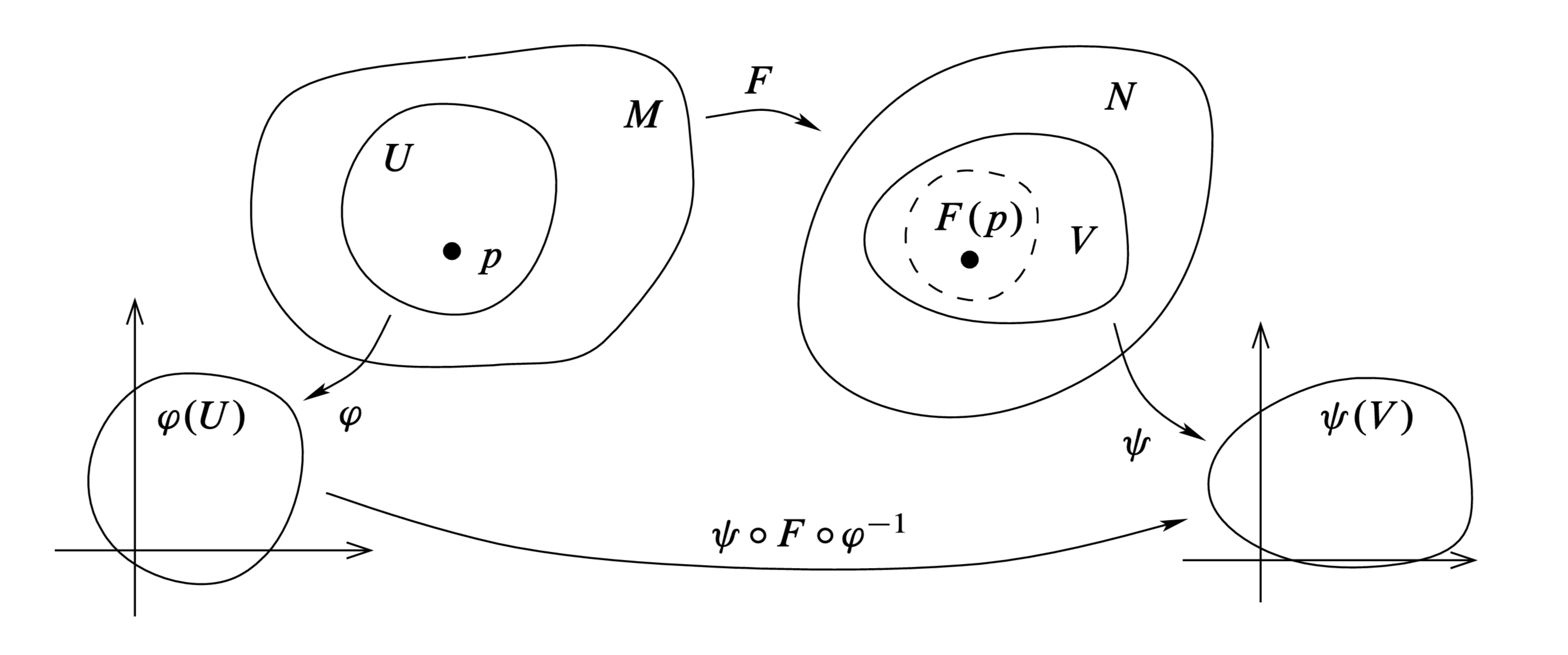}
\centering
\label{fig:mapeo-suave-en-cartas}
\end{figure}
\begin{definition}\label{notacion para mapeos de clase Ck}\index{notation for maps of class Ck@notation for maps of class \(C^k\)}
 If $M$ and $N$ are smooth manifolds with or without boundary and $k\in \mathbb{N}_{0}\cup\{\infty\}$, we write $C^{k}(M,N):=\{f\colon M\longrightarrow N\mid f\text{ is of class $C^{k}$}\}$. If $N=\mathbb{R}$, we simply write $C^{k}(M):=C^{k}(M,\mathbb{R})$.
\end{definition}
 \begin{remark}\label{obs:nociones-fundamentales-funcion-conoce-representacion-coordenada-ocasiones-denota}
 The function $\psi\circ F\circ\phi^{-1}$ is called the \textit{coordinate representation of $F$} and is sometimes denoted by $\hat{F}:=\psi\circ F\circ\phi^{-1}$.
 \end{remark}

\subsection{Scalar functions: values, charts, and regularity}

A scalar function on a manifold $M$ is, first of all, an assignment
$f\colon M\longrightarrow\mathbb R$: each point $p$ is assigned a unique
number $f(p)$. No coordinates, metric, or ambient space need to be chosen
for this assignment to make sense. On a surface, this is the scalar field
of vector calculus; examples include temperature or concentration
distributed over the surface. The manifold describes the domain, whereas
the function describes the datum that varies over it. Continuity and
differentiability are additional conditions on this assignment.

If $(U,\phi)$ is a smooth chart, its local expression is
\[
 f_\phi:=f\circ\phi^{-1}\colon\phi(U)\longrightarrow\mathbb R.
\]
The point $p$ and its coordinates $\phi(p)$ belong to different spaces,
but $f(p)=f_\phi(\phi(p))$ is the same value. For two charts
$(U,\phi)$ and $(V,\psi)$, one has
\begin{equation}\label{eq:compatibilidad-funciones-escalares-cartas}
 f_\psi=f_\phi\circ(\phi\circ\psi^{-1})
 \quad\text{on }\psi(U\cap V).
\end{equation}
Compatibility does not require identical formulas: it requires that
they assign the same number to the same point after changing its coordinates.

\begin{proposition}[Constructing a function from local expressions]
\label{prop:pegar-funciones-escalares-cartas}
Let $(U_\alpha,\phi_\alpha)_{\alpha\in A}$ be an atlas of $M$.
A family of functions
$f_\alpha\colon\phi_\alpha(U_\alpha)\longrightarrow\mathbb R$
determines a unique function $f\colon M\longrightarrow\mathbb R$ with
$f\circ\phi_\alpha^{-1}=f_\alpha$ if and only if
\[
 f_\beta(y)
 =f_\alpha\bigl(\phi_\alpha\circ\phi_\beta^{-1}(y)\bigr)
 \quad
 \bigl(y\in\phi_\beta(U_\alpha\cap U_\beta)\bigr)
\]
for all indices $\alpha,\beta$.
\end{proposition}
\begin{proof}
Necessity follows by substituting $y=\phi_\beta(p)$ into
$f_\alpha(\phi_\alpha(p))=f(p)=f_\beta(\phi_\beta(p))$.
Conversely, given $p\in M$, choose $\alpha$ with $p\in U_\alpha$ and
define $f(p):=f_\alpha(\phi_\alpha(p))$. If also $p\in U_\beta$,
compatibility gives the same value. Thus the assignment is well defined.
Its expression in each chart is the prescribed one, and any other
function with these expressions has the same value at every point.
\end{proof}

\begin{proposition}[Local criterion for scalar continuity and regularity]
\label{prop:regularidad-escalar-en-cartas}
Let $M$ be a smooth manifold with or without boundary,
$f\colon M\to\mathbb R$ a function, and $k\in\mathbb N_0\cup\{\infty\}$.
Then $f$ is continuous if and only if $f\circ\phi_\alpha^{-1}$ is continuous
in every chart of an atlas covering $M$.
Likewise, $f$ is of class $C^k$ if and only if these representations are
of class $C^k$. Any smooth atlas can be used; in particular, the condition
is equivalent to requiring it in every smooth chart.
For $k=0$, $C^0$ means continuity. In boundary charts,
for $k\geq1$ regularity is understood through local extension to Euclidean
space, in accordance with the convention already established.
\end{proposition}
\begin{proof}
If $f$ is continuous, its composition with the homeomorphism
$\phi_\alpha^{-1}$ is continuous. Conversely,
$f\restriction_{U_\alpha}=f_\alpha\circ\phi_\alpha$ is continuous.
For every open set $O\subseteq\mathbb R$,
\[
 f^{-1}(O)
 =\bigcup_{\alpha\in A}
 \bigl(f\restriction_{U_\alpha}\bigr)^{-1}(O)
\]
is open in $M$, since each $U_\alpha$ is open. This proves global continuity.

For $k\geq1$, the definition of regularity for maps, using the identity
chart of $\mathbb R$, reduces precisely to regularity of the scalar
expressions. Upon changing the chart,
identity~\eqref{eq:compatibilidad-funciones-escalares-cartas} expresses the
new representation as a composition with a smooth coordinate change.
The chain rule preserves the class $C^k$. At a boundary point, compose
local extensions, shrinking their domains if necessary; for $k=0$,
composition of continuous maps suffices.
Since the atlas domains cover $M$, the condition is local at every point.
For $k=\infty$, apply the same argument at each order.
\end{proof}

There are different ways to specify the assignment $f$. Given a smooth
realization $\iota\colon M\to\mathbb R^N$ and a function $F$ defined on a
neighborhood of $\iota(M)$, one can set $f=F\circ\iota$.
This is an extrinsic description; if $F$ is continuous or of class $C^k$,
the composition has the same regularity. The $N$ component functions of
$\iota$ may be globally defined without constituting a global chart of
the $n$-dimensional manifold: a chart must be a homeomorphism onto an open
subset of $\mathbb R^n$ or of the corresponding half-space.

One can also specify $f$ through compatible local expressions, without
choosing an embedding. Finally, a geometric property may directly produce
its values: after choosing a Riemannian metric and a point $q$, distance
to $q$ is an example studied in the following chapter.
Choosing intrinsic geometric data is not the same as choosing an ambient
space. These descriptions can represent the same function; the essential
object is the global assignment and its regularity.

\begin{example}[Height on the sphere and its stereographic expressions]
\label{ej:altura-esfera-cartas-autocontenido}
Consider
\[
 S^2:=\{(x,y,z)\in\mathbb R^3\mid x^2+y^2+z^2=1\},
 \qquad N:=(0,0,1),\quad S:=(0,0,-1).
\]
With the subspace topology, $S^2$ is Hausdorff and second countable.
Define $U_N:=S^2\setminus\{N\}$, $U_S:=S^2\setminus\{S\}$, and
\[
 \phi_N(x,y,z)=\left(\frac{x}{1-z},\frac{y}{1-z}\right),
 \qquad
 \phi_S(x,y,z)=\left(\frac{x}{1+z},\frac{y}{1+z}\right).
\]
Writing $r^2:=u^2+v^2$, their inverses are
\[
 \Phi_N(u,v)
 =\left(\frac{2u}{1+r^2},\frac{2v}{1+r^2},
              \frac{r^2-1}{1+r^2}\right),
 \qquad
 \Phi_S(u,v)
 =\left(\frac{2u}{1+r^2},\frac{2v}{1+r^2},
              \frac{1-r^2}{1+r^2}\right).
\]
Indeed, the sum of the squares of the three components is
$\displaystyle\frac{4r^2+(r^2-1)^2}{(1+r^2)^2}=1$.
Moreover, $1-z=\displaystyle\frac{2}{1+r^2}$ for $\Phi_N$ and
$1+z=\displaystyle\frac{2}{1+r^2}$ for $\Phi_S$, so that
$\phi_N\circ\Phi_N=\phi_S\circ\Phi_S=\operatorname{Id}_{\mathbb R^2}$.
For a point in $U_N$, the identity
$r^2=\displaystyle\frac{1+z}{1-z}$ recovers its three components through
$\Phi_N\circ\phi_N$; for $U_S$, use
$r^2=\displaystyle\frac{1-z}{1+z}$. All these expressions are continuous,
so both charts are homeomorphisms onto $\mathbb R^2$.

On the overlap, $r>0$ and
\[
 T(u,v):=\phi_S\circ\Phi_N(u,v)
 =\left(\frac{u}{r^2},\frac{v}{r^2}\right).
\]
This map is smooth on $\mathbb R^2\setminus\{0\}$, and $T\circ T$ is the
identity, since $|T(u,v)|^2=r^{-2}$. Thus the two charts determine a
smooth atlas covering the entire sphere.

The height function is $h(x,y,z):=z$. Its expressions are
\[
 h_N(u,v)=\frac{r^2-1}{1+r^2},
 \qquad
 h_S(u,v)=\frac{1-r^2}{1+r^2}.
\]
Compatibility is checked directly:
\[
 h_S(T(u,v))
 =\frac{1-r^{-2}}{1+r^{-2}}
 =\frac{r^2-1}{1+r^2}
 =h_N(u,v).
\]
Both formulas are smooth, including the origin in each chart. Hence $h$
is smooth also at the poles, where $h(N)=1$ and $h(S)=-1$.
For example,
\[
 \partial_u h_N=\frac{4u}{(1+r^2)^2},
 \qquad \partial_v h_N=\frac{4v}{(1+r^2)^2};
\]
the other chart gives the negatives of these expressions.

By contrast, writing the same formula $u$ in both charts would not
define a global function: for
$p=(\displaystyle\frac{4}{5},0,\displaystyle\frac{3}{5})$, one has
$\phi_N(p)=(2,0)$ and $\phi_S(p)=(\displaystyle\frac{1}{2},0)$.
The prescribed values would differ. Identical formulas alone do not
replace compatibility.

The word \emph{height} uses the inclusion in $\mathbb R^3$ and the chosen
axis. Changing the unit direction $a$ gives another function
$h_a(p):=\langle p,a\rangle$. Once chosen, each is a well-defined function
on $S^2$, although the abstract sphere singles out no axis.
Indeed, a function invariant under all rotations would be constant:
given $p,q\in S^2$, there is a rotation $R$ with $Rp=q$, and invariance
would give $f(q)=f(Rp)=f(p)$.
\end{example}

\begin{figura}[htbp]
 \includegraphics[width=0.88\linewidth]{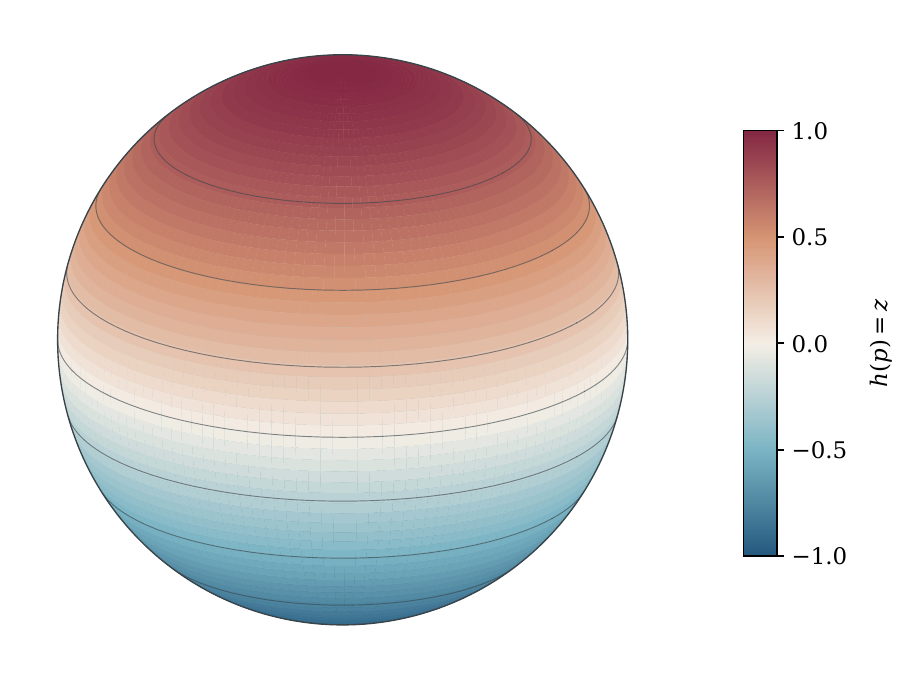}
 \caption{A scalar field on a surface: color represents the value
 $h(p)=z$, according to the indicated scale. Horizontal circles are
 level sets. The surface retains its shape; color represents a number
 assigned to each point, not an additional coordinate of the domain.}
 \label{fig:funcion-escalar-esfera-altura}
\end{figura}

\begin{example}[Continuity, failure of differentiability, and discontinuity]
\label{ej:regularidad-escalar-esfera}
On the same sphere, consider $a(p):=|h(p)|$ and
\[
 b(p):=
 \begin{cases}
  1,&h(p)>0,\\
  0,&h(p)\leq0.
 \end{cases}
\]
The function $a$ is continuous by composition. It is smooth where $h\neq0$,
but is not of class $C^1$ on the equator. To check this, fix
$p=(\cos\theta,\sin\theta,0)$ and use the smooth curve
\[
 c(t):=(\cos t\cos\theta,\cos t\sin\theta,\sin t).
\]
One has $c(0)=p$ and $(a\circ c)(t)=|\sin t|$. The quotients
$\displaystyle\frac{|\sin t|-0}{t}$ tend to $1$ from the right and to
$-1$ from the left. If $a$ were $C^1$ near $p$, its composition with $c$
would be $C^1$, contradicting these limits.

The function $b$ is smooth in each open hemisphere because it is constant
there, but it is discontinuous at every equatorial point:
$c(t)\to p$ as $t\to0$, whereas $b(c(t))=1$ for small $t>0$ and
$b(c(t))=0$ for small $t<0$. Changing the value of $b$ only on the equator
does not remove these two different limits. These are properties of the
functions, rather than defects of the charts used to represent them.
\end{example}

\subsection{Partitions of unity}

 A partition of unity turns local constructions into global objects while retaining control over supports. It will be our basic tool for gluing functions, metrics, sections, and estimates defined in different charts.

 \begin{definition}\label{def:nociones-fundamentales-espacio-topologico-cubierta-abierta-particion-unidad}\index{subordinate partition of unity@subordinate partition of unity}
 Let $M$ be a topological space and let $\mathscr{U}=\{U_{\alpha}\mid\alpha\in J\}$ be an open cover of $M$. A partition of unity subordinate to $\mathscr{U}$ is a family of continuous functions $(\psi_{\alpha})_{\alpha\in J}$, $\psi_{\alpha}\colon M\longrightarrow \mathbb{R}$, such that:
 \begin{enumerate}[label=(\alph*)]
 \item $0\leq \psi_{\alpha}(x)\leq 1$ for every $\alpha\in J$ and every $x\in M$.
 \item $\supp(\psi_{\alpha})\subseteq U_{\alpha}$ for each $\alpha\in J$.
 \item The family $(\supp(\psi_{\alpha}))_{\alpha\in J}$ is locally finite.
 \item $\displaystyle\sum_{\alpha\in J}\psi_{\alpha}(x)=1$ for each $x\in M$.
 \end{enumerate}
 \end{definition}

 We have the following well-known result.

 \begin{theorem}[Existence of partitions of unity]
 \label{particionesdelaunidad}\index{existence of partitions of unity}
 Suppose that $M$ is a smooth manifold with or without boundary and that $\mathscr{U}=(U_{\alpha})_{\alpha\in J}$ is an open cover of $M$. Then there exists a partition of unity consisting of smooth functions and subordinate to $\mathscr{U}$.
 \end{theorem}
\begin{proof}
We apply the construction in
Theorem~\ref{teo: variedad es paracompacta}, choosing, in each of the
compact sets occurring in its proof, two concentric regular coordinate
balls
\[
 B_i\Subset B_i'\Subset U_{a(i)},
 \qquad a(i)\in J,
\]
so that the inner balls $(B_i)_{i\geq1}$ cover $M$ and the family
of outer balls $(B_i')_{i\geq1}$ is locally finite. To do this,
for each point of the corresponding compact set we choose a
ball whose closure is contained in another ball in the same coordinate
domain, and then take a finite subcover. Since the outer
balls remain within the open sets $W_j$ used in that
proof, local finiteness is preserved.

Recall the elementary construction of a cutoff function. Let
\[
 \vartheta(t):=
 \begin{cases}
  e^{-\frac{1}{t}},&t>0,\\
  0,&t\leq0.
 \end{cases}
\]
All its derivatives at the origin vanish, so $\vartheta$ is
smooth. If $0<r<R_0<R$, the function
\[
 \chi_{r,R_0}(x):=
 \frac{\vartheta(R_0^2-\lVert x\rVert^2)}
 {\vartheta(R_0^2-\lVert x\rVert^2)
  +\vartheta(\lVert x\rVert^2-r^2)}
\]
is smooth, equals one on $\overline B_{\mathrm{euc}}(0,r)$, and has support
contained in
$\overline B_{\mathrm{euc}}(0,R_0)\subset B_{\mathrm{euc}}(0,R)$.
The denominator never vanishes. After translating the center in
coordinates, restricting this function to $\mathbb H^n$ gives the same
construction for a relative ball in a boundary chart.

Transporting these functions through the charts and extending them by zero,
we obtain $f_i\in C^\infty(M)$ such that
\[
 0\leq f_i\leq1,\qquad f_i=1\ \text{on }\overline B_i,\qquad
 \operatorname{supp}f_i\subseteq B_i'.
\]
The extension is smooth because the support of the coordinate expression is
contained in the interior of $B_i'$. The family of supports is locally
finite, so
\[
 f:=\displaystyle\sum_{i=1}^{\infty}f_i
\]
is a smooth function. Since the sets $B_i$ cover $M$, we have $f>0$. Hence
$g_i:=\frac{f_i}{f}$ is smooth, the family $(\operatorname{supp}g_i)$ is
locally finite, and $\displaystyle \sum_{i=1}^{\infty} g_i=1$.

For $\alpha\in J$, define
\[
 \psi_\alpha:=\displaystyle\sum_{\substack{i\in\mathbb N\\a(i)=\alpha}}g_i,
\]
interpreting the sum as zero if there are no corresponding indices. All
these sums are locally finite. Moreover, the union of the supports
appearing in each $\psi_\alpha$ is closed by local finiteness and is
contained in $U_\alpha$; consequently,
$\operatorname{supp}\psi_\alpha\subseteq U_\alpha$. It is immediate that
$0\leq\psi_\alpha\leq1$ and that
$\displaystyle \sum_{\alpha\in J}\psi_\alpha=1$. We have constructed the subordinate
partition. Restricting the Euclidean cutoff functions directly proves
the case with boundary.
\end{proof}

\begin{proposition}[Partitions of unity on manifolds with corners]
\label{prop:esquinas-particiones-unidad}
Every smooth manifold with corners is locally compact and paracompact,
admits a compact exhaustion
$K_j\subseteq\operatorname{Int}_M K_{j+1}$, and has smooth partitions
of unity subordinate to any open cover. If $K\subseteq U\subseteq M$,
with $K$ compact and $U$ open, there is $\chi\in C_c^\infty(U)$ such
that $0\leq\chi\leq1$ and $\chi=1$ on a neighborhood of $K$.
Interiors and supports are taken in $M$.
\end{proposition}

\begin{proof}
In a chart with image $\Omega\subseteq G_{n,r}$, choose a closed
Euclidean ball centered at the image of the point, small enough that
its intersection with $G_{n,r}$ is contained in $\Omega$. This
intersection is compact and contains a relative neighborhood of the
point. Transporting it proves local compactness of $M$. Second
countability gives a countable basis of precompact open sets $(V_j)$.
Start with $K_1=\overline V_1$; to pass from $K_j$ to $K_{j+1}$, cover
$K_j\cup\overline V_{j+1}$ by finitely many precompact open sets and
take the union of their closures. Then
$K_j\cup\overline V_{j+1}\subseteq\operatorname{Int}_M K_{j+1}$ and
$M=\bigcup_j K_j$.

Set $K_0=\varnothing$,
$C_0=K_1$, $W_0=\operatorname{Int}_M K_2$, and, for $j\geq1$,
\[
 C_j:=K_{j+1}\setminus\operatorname{Int}_M K_j,
 \qquad W_j:=\operatorname{Int}_M K_{j+2}\setminus K_{j-1}.
\]
The compact sets $C_j$ cover $M$ and lie in the open sets $W_j$.
At each point of $C_j$, choose two concentric relative coordinate
balls $B_i\Subset B_i'$, with the closure of the larger one contained
in $W_j$ and in a member of the given cover. A finite subfamily of
the smaller balls covers $C_j$. Combining these families gives a
countable cover; the larger balls are locally finite. Indeed, if a
point lies in $\operatorname{Int}_M K_q$, this neighborhood meets no
$W_j$ with $j\geq q+1$, and only finitely many balls occur in the
remaining indices. This also proves paracompactness.

The radial cutoff in the proof of Theorem~\ref{particionesdelaunidad},
translated to each center and restricted to $G_{n,r}$, is one on the
smaller ball and has compact support within the larger one. Transporting
it and extending by zero gives $f_i\in C_c^\infty(M)$ with
$0\leq f_i\leq1$, $f_i=1$ on $B_i$, and
$\operatorname{supp}f_i\subseteq B_i'$. Extension by zero is smooth
because the support lies within the relatively open chart domain.
The locally finite sum $f=\sum_i f_i$ is smooth and positive.
The functions $g_i=f_i/f$ sum to one. Grouping them by the cover member
containing $B_i'$ gives a subordinate partition: the union of supports
in each group is closed by local finiteness and stays inside that
member. Hence the support of the sum also stays there.

For the last assertion, cover $K$ by finitely many smaller balls whose
corresponding larger balls have compact closure in $U$, and take the
preceding cutoffs $f_1,\ldots,f_s$. The function
$\chi=1-\prod_{i=1}^s(1-f_i)$ is smooth, takes values between zero
and one, equals one on the union of the smaller balls, and has support
in the finite union of the compact supports, contained in $U$.
\end{proof}

\begin{proposition}[Regular charts adapted to depth]
\label{prop:esquinas-cartas-regulares-cubierta}
Let $M$ be a smooth manifold with corners of dimension $n$. Every point
$p\in M$, of depth $d$, has a chart $(V,\phi)$ such that
\[
 \phi(p)=0,\qquad
 \phi(V)=B_{\mathrm{euc}}(0,2R)\cap G_{n,d}
\]
for some $R>0$. The open set
$U:=\phi^{-1}(B_{\mathrm{euc}}(0,R)\cap G_{n,d})$ has compact closure
contained in $V$; we call $(U,\phi|_U)$ a regular chart with corners.
These charts form a basis for the topology.

Given an open cover by trivialization domains of a vector bundle,
there is a countable cover $(U_\alpha)_{\alpha\in A}$ by regular charts,
with corresponding larger open sets $V_\alpha$, such that
$(V_\alpha)_{\alpha\in A}$ is locally finite and each $V_\alpha$
is contained in one of those domains. There is a smooth partition
of unity $(\psi_\alpha)$ with
$\operatorname{supp}\psi_\alpha\Subset U_\alpha$.
When $M$ is compact, finitely many indices suffice.
\end{proposition}

\begin{proof}
Choose an initial chart with image $\Omega\subseteq G_{n,r}$ and
let $z$ be the image of $p$. Of the $r$ nonnegative coordinates,
exactly $d$ vanish at $z$. Shrink the domain so that the remaining
$r-d$ coordinates stay positive. A permutation places the vanishing
coordinates last, and translation of the remaining $n-d$ coordinates
maps $z$ to the origin. If the radius is smaller than half of each
initially positive coordinate, its positivity restriction does not
affect the ball. Thus the model in this ball is exactly $G_{n,d}$.
Relative openness allows us to choose $R>0$ so that
\[
 \overline B_{\mathrm{euc}}(0,3R)\cap G_{n,d}
\]
lies in the reduced image. Restrict the chart to the ball of radius
$2R$. The closure of $U$ is the inverse image of the compact set
$\overline B_{\mathrm{euc}}(0,R)\cap G_{n,d}$ and is contained in $V$.
The permutations and translations used have smooth affine extensions
with smooth inverses, so they preserve the structure. We can shrink
$R$ so that $V$ lies in any previously chosen neighborhood of $p$.
This proves the basis assertion.

To construct the cover, use the exhaustion and the compact sets $C_j$
contained in the open sets $W_j$ from the proof of
Proposition~\ref{prop:esquinas-particiones-unidad}. For every
$p\in C_j$, choose a pair $U\Subset V$ of the preceding form, with
$\overline V$ contained in $W_j$ and in a trivialization domain.
Compactness of $C_j$ gives a finite family of the smaller sets
covering it. Their union is countable and covers $M$.
If $p\in\operatorname{Int}_M K_q$, this neighborhood does not meet
$W_j$ for $j\geq q+1$, and only finitely many pairs remain at the
earlier indices. Therefore the family of larger sets is locally
finite. Apply the cited proposition to the cover $(U_\alpha)$.
The support of $\psi_\alpha$ is closed in $M$, lies in $U_\alpha$,
and is compact because it is a closed subset of
$\overline{U_\alpha}$. Hence
$\operatorname{supp}\psi_\alpha\Subset U_\alpha$.
If $M$ is compact, a finite subcover by regular charts and a
subordinate partition give the last assertion.
\end{proof}

 The existence of partitions of unity has the following important consequence:
 \begin{proposition}[Existence of bump functions]\label{funciones flan}\index{existence of bump functions}
 Let $M$ be a smooth manifold with or without boundary. For any closed set $A\subseteq M$ and any open set $U\subseteq M$ such that $A\subseteq U$, there exists a smooth function $\psi\colon M\longrightarrow \mathbb{R}$ such that $0\leq \psi\leq 1$, $\psi\equiv 1$ on $A$, and $\supp(\psi)\subseteq U$.

 \end{proposition}
\begin{proof}
Apply Theorem~\ref{particionesdelaunidad} to the cover
$\{U,M\setminus A\}$. Let $(\psi_0,\psi_1)$ be a subordinate partition,
with $\operatorname{supp}\psi_0\subseteq U$ and
$\operatorname{supp}\psi_1\subseteq M\setminus A$. On $A$ we have
$\psi_1=0$ and, since $\psi_0+\psi_1=1$, we have $\psi_0=1$. Therefore,
$\psi:=\psi_0$ has all the required properties.
\end{proof}
 This result allows us to extend smooth functions defined in a neighborhood of a given point to all of $M$ as follows:
 \begin{theorem}[Extension of smooth functions]\label{extension de funciones suaves}\index{extension of smooth functions@extension of smooth functions}
 Let $M$ be a smooth manifold with or without boundary and let $U\subseteq M$ be open. If $p\in U$ and $f\in C^{\infty}(U)$, then there exist an open set $W\subseteq M$ such that $p\in W\subseteq\overline{W}\subseteq U$ and a smooth function $\tilde{f}\colon M\longrightarrow \mathbb{R}$ such that $\tilde{f}\restriction_{W}=f$.
 \end{theorem}
 \begin{proof}
 Let $U\subseteq M$ be open, and let $p\in U$ and $f\in C^{\infty}(U)$. Since $M$ is a regular topological space, there exists an open set $W\subseteq M$ such that $p\in W\subseteq \overline{W}\subseteq U$. By Proposition~\ref{funciones flan}, there exists a smooth function $\psi\colon M\longrightarrow \mathbb{R}$ such that $0\leq \psi\leq 1$, $\psi\equiv 1$ on $\overline{W}$, and $\supp(\psi)\subseteq U$. Define \[\tilde{f}\colon M\longrightarrow \mathbb{R}, \qquad \tilde{f}(q)=\begin{cases}
 \psi(q)f(q)&\text{if $q\in U$}\\
 0&\text{if $q\notin U$}
 \end{cases}.\] Then $\tilde{f}$ is smooth because $\supp(\psi )\subseteq U$ and $\tilde{f}\restriction_{W}=f$.
 \end{proof}
 \section{The tangent bundle}

Approximating a nonlinear object by linear data is one of the central ideas of analysis. In Example~\ref{ejemplo plano tangente}, we saw that the graph of a differentiable function is approximated near each point by its tangent plane. On an abstract manifold, we need to preserve this intuition without relying on an ambient space.

 In multivariable calculus, many manifolds can be described locally as the image of a \textit{parametrization}
\[
\psi\colon \widetilde{U} \subset \mathbb{R}^{k} \longrightarrow \mathbb{R}^{n},
\]
where $\widetilde{U}\subseteq \mathbb{R}^{k}$ is open, $k \leq n$, and $d\psi_{x}$ has rank $k$ for each $x\in \widetilde{U}$. For example, the sphere, cylinder, and torus admit descriptions of this kind in small neighborhoods of each point.

Fix $x_0 \in \widetilde U$ and let $p = \psi(x_0)$. The Fréchet derivative of $\psi$ at $x_0$,
\[
\mathrm{d}\psi_{x_0}\colon \mathbb{R}^{k} \longrightarrow \mathbb{R}^{n},
\]
is a linear map that describes how the coordinates of $\mathbb{R}^{n}$ change when we vary the coordinates of $\widetilde U$ slightly around $x_0$. The image of $\mathrm{d}\psi_{x_0}$ is a vector subspace of $\mathbb{R}^{n}$ of dimension $k$ that ``best approximates'' $\psi(\widetilde{U})$ near $p$. The following example illustrates this:

\begin{example}\label{ej:nociones-fundamentales-parametrizacion-funcion-describe-cilindro-radio-fijemos}
Consider the parametrization $\psi\colon \mathbb{R}^{2} \longrightarrow \mathbb{R}^{3}$ given by
\[
\psi(x,y) = (\cos x,\ \sin x,\ y).
\]
This function describes a cylinder of radius $1$ in $\mathbb{R}^{3}$. Fix ${x_0}=(0,0)$, so that $p = (1,0,0)$. The Fréchet derivative at $(0,0)$ is
\[
\mathrm{d}\psi_{(0,0)}(h_{1},h_{2}) = (-\sin 0\cdot h_{1},\ \cos 0\cdot h_{1},\ h_{2}) = (0,\ h_{1},\ h_{2}).
\]
The space \[\text{Im}(\mathrm{d}\psi_{(0,0)}):=\{(0,y,z)\mid y,z\in \mathbb{R}\}\] is a vector subspace of $\mathbb{R}^{3}$ of dimension $2$. Translating it by the vector $(1,0,0)$ gives the tangent plane to the cylinder at the point $(1,0,0)$, given by $\Pi:=\{(1,y,z)\in \mathbb{R}^{3}\mid y,z\in \mathbb{R}\}$, as illustrated below:
\begin{center}
 \includegraphics[scale=1]{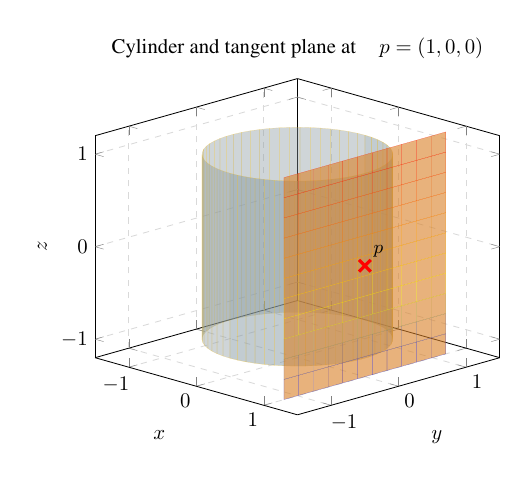}
\end{center}
This subspace has basis \[\beta=\displaystyle\left\{\frac{\partial \psi}{\partial x}(x_{0}),\frac{\partial\psi}{\partial y}(x_{0})\right\}=\left\{(0,1,0),(0,0,1)\right\}.\]
\end{example}
Under these assumptions, we refer to
$\operatorname{Im}(d\psi_{x_0})$ as the tangent space to
the local image of $\psi$ around $x_0$, at the point $p=\psi(x_0)$.

The local tool that makes this description precise is the following
theorem from differential calculus.
\begin{theorem}[Rank theorem]\label{teo: del rango euclidiano}\index{rank theorem}
Let \(f\colon \widetilde U\subseteq \mathbb{R}^k\longrightarrow \mathbb{R}^n\) be a function of class \(C^1\), and let \(x_0\in \widetilde U\). Suppose there exists \(r>0\), with \(B_{\mathrm{euc}}(x_0,r)\subseteq \widetilde U\), such that \(df_x\) has constant rank \(m\) for every \(x\in B_{\mathrm{euc}}(x_0,r)\). Then there exist open sets \(\widetilde V_1,\widetilde V_2\subseteq \mathbb{R}^k\) and \(V_1,V_2\subseteq \mathbb{R}^n\), with \(x_0\in \widetilde V_1\), \(0\in \widetilde V_2\), \(f(x_0)\in V_1\), and \(0\in V_2\), and diffeomorphisms of class \(C^1\)
\[
\Phi\colon \widetilde V_2\longrightarrow \widetilde V_1,
\qquad
\Psi\colon V_2\longrightarrow V_1,
\]
with $\widetilde V_1\subseteq\widetilde U$,
$f(\widetilde V_1)\subseteq V_1$, $\Phi(0)=x_0$, and $\Psi(0)=f(x_0)$,
such that
\[
\Psi^{-1}\circ f\circ \Phi(x_1,\dots,x_k)
=
(x_1,\dots,x_m,0,\dots,0).
\]
In particular, near \(x_0\), the map \(f\) can be expressed in suitable coordinates on the domain and codomain as projection onto the first \(m\) coordinates followed by the natural inclusion into \(\mathbb{R}^n\).
If $f$ is of class $C^j$, where $1\leq j\leq\infty$, then
$\Phi$ and $\Psi$ may be chosen of class $C^j$.
\end{theorem}
\begin{proof}
If $m=0$, then $df_x=0$ on $B_{\mathrm{euc}}(x_0,r)$. Integrating along
the segments joining $x_0$ to points of this ball gives $f(x)=f(x_0)$.
It suffices to take $\widetilde V_2=B_{\mathrm{euc}}(0,r)$,
$\widetilde V_1=B_{\mathrm{euc}}(x_0,r)$, $V_2=\mathbb R^n$,
$V_1=\mathbb R^n$, and the translations
$\Phi(y)=x_0+y$, $\Psi(z)=f(x_0)+z$.

Assume $m\geq1$. After permuting coordinates in the domain and codomain,
we may assume that the submatrix formed by the first $m$ rows and columns
of $df_{x_0}$ is invertible. Write $f=(f^1,\ldots,f^n)$ and define
\[
H(x):=\bigl(f^1(x)-f^1(x_0),\ldots,f^m(x)-f^m(x_0),
             x^{m+1}-x_0^{m+1},\ldots,x^k-x_0^k\bigr).
\]
Then $H(x_0)=0$ and
\[
dH_{x_0}=\begin{pmatrix}A&B\\0&I_{k-m}\end{pmatrix},
\]
where $A$ is the chosen invertible submatrix.
Theorem~\ref{lem:funcion-inversa-euclidea-previa-variedades} provides a
local inverse $\Phi$ of $H$, of class $C^1$, with $\Phi(0)=x_0$.
Restrict its domain to a product of open boxes
$\widetilde V_2=P\times Q\subseteq\mathbb R^m\times\mathbb R^{k-m}$
containing the origin and whose image
$\widetilde V_1:=\Phi(P\times Q)$ lies in
$B_{\mathrm{euc}}(x_0,r)$. If $m=k$, the second factor is $\mathbb R^0$.

For $y=(y',y'')\in P\times Q$, the identity $H\circ\Phi(y)=y$ gives
\[
f\circ\Phi(y)=\bigl(y'+c,h(y',y'')\bigr),
\qquad
c:=\bigl(f^1(x_0),\ldots,f^m(x_0)\bigr),
\]
where $h$ takes values in $\mathbb R^{n-m}$. By the chain rule,
$d(f\circ\Phi)_y$ has rank $m$, and its matrix is
\[
\begin{pmatrix}I_m&0\\D_1h(y)&D_2h(y)\end{pmatrix}.
\]
The first $m$ rows are linearly independent. The remaining rows are
linear combinations of them, so $D_2h(y)=0$.
Since $Q$ is convex and contains the origin,
\[
h(y',y'')-h(y',0)=\int_0^1D_2h(y',ty'')y''\,dt=0.
\]
Define $b(y'):=h(y',0)$. Then
\[
f\circ\Phi(y',y'')=\bigl(y'+c,b(y')\bigr).
\]

Let $V_2:=P\times\mathbb R^{n-m}$ and
$V_1:=(c+P)\times\mathbb R^{n-m}$. The map
\[
\Psi\colon V_2\longrightarrow V_1,
\qquad
\Psi(z',z''):=\bigl(c+z',z''+b(z')\bigr)
\]
is a diffeomorphism, since its inverse is given explicitly by
\[
\Psi^{-1}(w',w'')=\bigl(w'-c,w''-b(w'-c)\bigr).
\]
Moreover, $\Psi(0)=f(x_0)$, $f(\widetilde V_1)\subseteq V_1$, and
\[
\Psi^{-1}\circ f\circ\Phi(y',y'')=(y',0),
\]
as required. The initial permutations are incorporated into $\Phi$ and
$\Psi$ when returning to the original coordinates. If $f\in C^j$, the
inverse function theorem gives $\Phi\in C^j$; composition and restriction
give $b\in C^j$, and the formulas for $\Psi$ and $\Psi^{-1}$ give the
same regularity. This includes $j=\infty$.
\end{proof}
\begin{proposition}[Local description by curves]
\label{prop tangente y parametrizacion}
Let $\psi\colon\widetilde U\subseteq\mathbb R^k\to\mathbb R^n$ be of class
$C^1$, let $x_0\in\widetilde U$, and assume that
$d\psi_{x_0}$ has rank $k$. Then there exists an open neighborhood
$V\subseteq\widetilde U$ of $x_0$ such that
$\psi\restriction_V\colon V\to\mathbb R^n$ is an embedding of class $C^1$: it is an injective immersion and a
homeomorphism onto its image with the subspace topology.
If $p=\psi(x_0)$, every such neighborhood satisfies
\[
 \operatorname{Im}(d\psi_{x_0})
 =
 \left\{\gamma'(0)\ \middle|\
 \begin{array}{l}
 \gamma\colon(-\varepsilon,\varepsilon)\to\mathbb R^n
       \text{ is of class }C^1,\quad \gamma(0)=p,\\
 \gamma((-\varepsilon,\varepsilon))\subseteq\psi(V)
       \text{ for some }\varepsilon>0
 \end{array}\right\}.
\]
\end{proposition}
\begin{proof}
Since $d\psi_{x_0}$ has rank $k$, we can choose a coordinate projection
$\pi\colon\mathbb R^n\to\mathbb R^k$ such that
$d(\pi\circ\psi)_{x_0}$ is invertible.
Theorem~\ref{lem:funcion-inversa-euclidea-previa-variedades} provides open
sets $V\subseteq\widetilde U$ and $A\subseteq\mathbb R^k$, with
$x_0\in V$ and $\pi(p)\in A$, such that
$\pi\circ\psi\restriction_V\colon V\to A$ is a diffeomorphism of class
$C^1$. Denote its inverse by $\theta\colon A\to V$.
The identity
\[
\theta\circ\pi\circ\psi\restriction_V=\operatorname{Id}_V
\]
shows that $\psi\restriction_V$ is injective. Differentiating it shows
that $d\psi_x$ is injective for every $x\in V$. Moreover, the restriction
of $\theta\circ\pi$ to $\psi(V)$ is a continuous inverse of
$\psi\restriction_V$. Hence this restriction is an embedding of class $C^1$.

If $v=d\psi_{x_0}(w)$, choose $\varepsilon>0$ such that $x_0+tw\in V$
for $|t|<\varepsilon$. The curve $\gamma(t):=\psi(x_0+tw)$ satisfies
$\gamma(0)=p$, has image in $\psi(V)$, and satisfies $\gamma'(0)=v$.
The same formula gives the constant curve when $w=0$.

Conversely, let $\gamma$ be a curve as in the statement. The map
$\theta\circ\pi$ is defined and of class $C^1$ on the open set
$\pi^{-1}(A)\subseteq\mathbb R^n$. Thus
\[
\delta(t):=\theta(\pi(\gamma(t)))
\]
is a curve of class $C^1$ in $V$, with $\delta(0)=x_0$ and
$\psi\circ\delta=\gamma$. The chain rule gives
\[
\gamma'(0)=d\psi_{x_0}(\delta'(0))\in\operatorname{Im}(d\psi_{x_0}).
\]

Finally, let $V'$ be any other open neighborhood of $x_0$ on which $\psi$
is an embedding. The set $V\cap V'$ is open in $V'$ and contains $x_0$.
Since $\psi\restriction_{V'}$ is a homeomorphism onto its image,
$\psi(V\cap V')$ is a neighborhood of $p$ in $\psi(V')$.
Every continuous curve in $\psi(V')$ passing through $p$ remains in this
neighborhood on a sufficiently small interval. The preceding argument
applies to its restriction; the other inclusion follows from
$\psi(x_0+tw)$ with $x_0+tw\in V'$. Hence the equality also holds for $V'$.
\end{proof}

\begin{remark}[Why the assertion is local]
We cannot replace $\psi(V)$ by the entire image $\psi(\widetilde U)$.
An immersion may have self-intersections: a curve passing through
$p=\psi(x_0)$ along another branch may have a velocity that does not
belong to $\operatorname{Im}(d\psi_{x_0})$. The choice of $V$ separates
the branch determined by $x_0$.
\end{remark}
In what follows, we return to considering manifolds in $\mathbb{R}^{n}$; taking $k\leq n$ and considering it to have dimension $k$ served to illustrate what happens when the manifold lies in a possibly larger ambient space.

The preceding motivation leads to an intrinsic definition that works
unchanged at interior and boundary points alike.
\begin{definition}\label{def canonica espacio tangente derivaciones en p}\index{tangent space via derivations}
Let $M$ be a smooth manifold with or without boundary and let $p\in M$. A
\textbf{derivation at $p$} is a linear map
$v\colon C^{\infty}(M)\longrightarrow\mathbb R$ such that
\[
 v(fg)=f(p)v(g)+g(p)v(f)
\]
for any $f,g\in C^\infty(M)$. The set of all
derivations at $p$ is denoted by $T_pM$ and is called the
\textbf{tangent space to $M$ at $p$}. Its elements are the tangent
vectors at $p$.
\end{definition}
\begin{remark}\label{obs:nociones-fundamentales-espacio-tangente-espacio-vectorial-real-operaciones}
With addition and scalar multiplication defined pointwise, $T_pM$ is
a real vector space. If $v\in T_pM$, then $v(1)=0$ and therefore
$v(c)=0$ for every constant function $c$.
\end{remark}

At interior points, we recover the geometric description using
two-sided curves.
\begin{definition}[Tangent space defined by curves]\label{def:nociones-fundamentales-espacio-tangente-definido-por-curvas}\index{tangent space defined by curves}
Let $M$ be a smooth manifold with or without boundary and let $p\in\operatorname{Int}(M)$.
Consider smooth curves
$\gamma\colon(-\varepsilon,\varepsilon)\longrightarrow M$ with
$\gamma(0)=p$. Two such curves are equivalent if, for some interior chart
$\varphi$ around $p$, we have
$(\varphi\circ\gamma_1)'(0)=(\varphi\circ\gamma_2)'(0)$. The chain
rule shows that this condition is independent of the chart. We denote the
class of $\gamma$ by $[\gamma]$.
\end{definition}

\begin{lemma}[Locality of derivations]
\label{lem:localidad-derivaciones-punto}
Let $M$ be a smooth manifold with or without boundary, let $p\in M$, and let
$v\in T_pM$. If $f,g\in C^\infty(M)$ agree in a neighborhood of
$p$, then $v(f)=v(g)$. Consequently, $v$ can be evaluated unambiguously
on the germ of a smooth function defined only near $p$.
\end{lemma}
\begin{proof}
Set $h=f-g$ and choose an open set $W$ containing $p$ on which
$h=0$. By Proposition~\ref{funciones flan}, there exists
$\chi\in C^\infty(M)$ with $\chi(p)=1$ and
$\operatorname{supp}\chi\subseteq W$. Since $\chi h=0$ and $h(p)=0$, the
Leibniz rule gives
\[
 0=v(\chi h)=\chi(p)v(h)+h(p)v(\chi)=v(h).
\]
Therefore, $v(f)=v(g)$. If a function is defined only near
$p$, Theorem~\ref{extension de funciones suaves} provides a
global extension agreeing with it on a smaller neighborhood; the equality
just proved shows that the value is independent of the extension.
\end{proof}

The identification between curves and derivations uses the following two
formulas from Euclidean calculus.
\begin{theorem}[Taylor's theorem with multi-indices]\label{taylor multivariable multiindices}\index{Taylor's theorem@Taylor's theorem!with multi-indices}
 Let $f\in C^{k+1}(U)$ and $a\in U$. We define the Taylor polynomial of $f$ of order $k$ centered at $a\in U$ by
 \[
P_{k,a,f}(x):=\displaystyle\sum_{|\alpha|\leq k}\frac{1}{\alpha !}D^{\alpha}f(a)(x-a)^{\alpha},
 \]
 where $\displaystyle\frac{1}{\alpha !}:=\displaystyle\frac{1}{\alpha_1 !}\cdots \displaystyle\frac{1}{\alpha_n !}$. Define $R_{k,a,f}(x)=f(x)-P_{k,a,f}(x)$ and call it the remainder of order $k$ of $f$ at $a$. Then:
 \begin{enumerate}
 \item If $x\in U$ is such that $[a,x]\subseteq U$, there exists $\xi\in[a,x]$ such that
 \[
 R_{k,a,f}(x)=\displaystyle\sum_{|\alpha| = k+1} \frac{1}{\alpha!} D^{\alpha}f\big(\xi\big) (x - a)^{\alpha}.
 \]
 \item If $x\in U$ is such that $[a,x]\subseteq U$, then
 \[
 R_{k,a,f}(x)=\displaystyle\sum_{|\alpha|=k+1}\frac{k+1}{\alpha !}(x-a)^{\alpha}\int_{0}^{1}(1-t)^{k}D^{\alpha}f(a+t(x-a))dt.
 \]
 \item We have
 \[
 \lim_{x\to a}\frac{f(x)-P_{k,a,f}(x)}{\|x-a\|^{k}}
 =\lim_{x\to a}\frac{R_{k,a,f}(x)}{\|x-a\|^{k}}=0
 \]
 and $P_{k,a,f}(x)$ is the unique polynomial of degree at most $k$ with this property.
 \end{enumerate}
\end{theorem}
\begin{proof}
Fix \(x\in U\) with \([a,x]\subseteq U\), set \(h=x-a\), and
define \(\gamma(t):=f(a+th)\). Induction on \(j\), together with the chain
rule, the Clairaut--Schwarz theorem, and the multinomial formula, gives
\begin{equation}
\label{eq:derivadas-curva-taylor-multiindice}
\gamma^{(j)}(t)
=\displaystyle\sum_{|\alpha|=j}\frac{j!}{\alpha!}
D^\alpha f(a+th)h^\alpha,
\qquad 0\leq j\leq k+1.
\end{equation}
Indeed, if the equality holds for \(j\), differentiation gives
\[
\displaystyle\sum_{|\alpha|=j}\displaystyle\sum_{r=1}^n
\frac{j!}{\alpha!}D^{\alpha+e_r}f(a+th)h^{\alpha+e_r}.
\]
For a multi-index \(\beta\) of length \(j+1\), the coefficient of the
corresponding term is
\[
\displaystyle\sum_{r:\,\beta_r>0}\frac{j!}{(\beta-e_r)!}
=\frac{j!}{\beta!}\displaystyle\sum_{r=1}^n\beta_r
=\frac{(j+1)!}{\beta!},
\]
which completes the induction.

The one-dimensional Taylor formula with integral remainder, obtained by
successive applications of the fundamental theorem of calculus, states that
\[
\gamma(1)=\displaystyle\sum_{j=0}^k\frac{\gamma^{(j)}(0)}{j!}
+\frac1{k!}\int_0^1(1-t)^k\gamma^{(k+1)}(t)\,dt.
\]
Substituting \eqref{eq:derivadas-curva-taylor-multiindice}, the finite sum
is \(P_{k,a,f}(x)\), and the last term is
\[
\displaystyle\sum_{|\alpha|=k+1}\frac{k+1}{\alpha!}h^\alpha
\int_0^1(1-t)^kD^\alpha f(a+th)\,dt.
\]
This proves (2). The mean value theorem for integrals gives
\(\theta\in[0,1]\) such that the one-dimensional remainder is
\(\gamma^{(k+1)}(\theta)/(k+1)!\). Using
\eqref{eq:derivadas-curva-taylor-multiindice} again and taking
\(\xi=a+\theta h\) gives (1).

The derivatives of order \(k+1\) are bounded on a closed ball centered
at \(a\) and contained in \(U\). For \(x\)
near \(a\), the integral formula implies
\[
|R_{k,a,f}(x)|\leq C\|x-a\|^{k+1},
\]
and proves the limit in (3). If \(Q\) is another polynomial of degree at most
\(k\) with this property, then
\[
\frac{(Q-P_{k,a,f})(x)}{\|x-a\|^k}\longrightarrow0.
\]
If \(S:=Q-P_{k,a,f}\ne0\), write
\(S(a+h)=\sum_{\ell=0}^k S_\ell(h)\), where each \(S_\ell\) is homogeneous
of degree \(\ell\). Let \(S_j\) be the first nonzero part and choose
\(v\) with \(S_j(v)\ne0\). Then, for \(t>0\),
\[
\frac{S(a+tv)}{\|tv\|^k}
=\frac{t^{j-k}}{\|v\|^k}\left(S_j(v)+\sum_{\ell=j+1}^k t^{\ell-j}S_\ell(v)\right)
\]
does not tend to zero as \(t\to0^+\): it is unbounded if \(j<k\), and
tends to \(S_k(v)/\|v\|^k\ne0\) if \(j=k\). This contradiction proves
uniqueness.
\end{proof}

\begin{lemma}[Finite-dimensional Hadamard lemma with integral remainder]
\label{lem:hadamard-finito-dimensional}
\index{Hadamard lemma}
Let $U\subseteq\mathbb R^n$ be open, let $x_0,x\in U$, and suppose that
$x_0+t(x-x_0)\in U$ for every $t\in[0,1]$. If
$F\colon U\longrightarrow\mathbb R^N$ is of class $C^1$, then
\begin{equation}
\label{eq:hadamard-primer-orden-finito-dimensional}
F(x)-F(x_0)
=
\int_0^1
DF\bigl(x_0+t(x-x_0)\bigr)[x-x_0]\,dt.
\end{equation}
In particular,
\begin{equation}
\label{eq:hadamard-factorizacion-coordenada}
F(x)-F(x_0)
=
\displaystyle\sum_{j=1}^n(x^j-x_0^j)F_j(x,x_0),
\qquad
F_j(x,x_0)
:=
\int_0^1
\frac{\partial F}{\partial x^j}
\bigl(x_0+t(x-x_0)\bigr)\,dt.
\end{equation}
If $F$ is of class $C^2$, $F(x_0)=0$, and $DF(x_0)=0$, then
\begin{equation}
\label{eq:hadamard-segundo-orden-finito-dimensional}
F(x)
=
\displaystyle\sum_{j,k=1}^n
(x^j-x_0^j)(x^k-x_0^k)
\int_0^1(1-t)
\frac{\partial^2F}{\partial x^j\partial x^k}
\bigl(x_0+t(x-x_0)\bigr)\,dt.
\end{equation}
The same formulas hold for a family $F(z,x)$ depending smoothly
on a finite-dimensional parameter $z$: the functions defined by
the integrals above are jointly smooth in $(z,x,x_0)$ as long as
the segments under consideration remain in the domain.
\end{lemma}

\begin{proof}
Define $g(t):=F(x_0+t(x-x_0))$. The chain rule gives
\[
g'(t)=DF\bigl(x_0+t(x-x_0)\bigr)[x-x_0].
\]
The fundamental theorem of calculus, applied componentwise in
$\mathbb R^N$, proves
$\displaystyle F(x)-F(x_0)=g(1)-g(0)=\int_0^1g'(t)\,dt$ and, upon expressing the differential in
the standard basis, yields
\eqref{eq:hadamard-factorizacion-coordenada}.

Now suppose that $F(x_0)=0$ and $DF(x_0)=0$. For each $t\in[0,1]$, applying the
preceding formula to $DF$ on the subsegment joining $x_0$ to
$x_0+t(x-x_0)$ gives
\[
DF\bigl(x_0+t(x-x_0)\bigr)[x-x_0]
=
t\int_0^1
D^2F\bigl(x_0+st(x-x_0)\bigr)[x-x_0,x-x_0]\,ds.
\]
Substitute this identity into
\eqref{eq:hadamard-primer-orden-finito-dimensional}. The change of variable
$\tau=st$ and Fubini's theorem on the triangle
$0\leq\tau\leq t\leq1$ yield
\[
F(x)
=
\int_0^1(1-\tau)
D^2F\bigl(x_0+\tau(x-x_0)\bigr)[x-x_0,x-x_0]\,d\tau,
\]
which is precisely
\eqref{eq:hadamard-segundo-orden-finito-dimensional} after expanding the
bilinear form in coordinates. The assertion with parameters follows
by differentiating under the integral over the compact interval $[0,1]$.
\end{proof}

\begin{proposition}[Equivalence between curves and derivations]
\label{prop:equivalencia-curvas-derivaciones-tangentes}
Let \(M\) be a smooth manifold with or without boundary and let
\(p\in\operatorname{Int}(M)\). The map
\[
\{\text{classes of two-sided curves at \(p\)}\}
\longrightarrow T_pM,
\qquad
[\gamma]\longmapsto v_{[\gamma]},
\]
where
\[
v_{[\gamma]}(f):=
\left.\frac d{dt}(f\circ\gamma)(t)\right|_{t=0},
\qquad f\in C^\infty(M),
\]
is a canonical bijection.
\end{proposition}
\begin{proof}
If \([\gamma_1]=[\gamma_2]\), choose an interior chart
\(\varphi\colon U\to\mathbb R^n\) around \(p\). For
\(F=f\circ\varphi^{-1}\), the chain rule gives
\[
\left.\frac d{dt}(f\circ\gamma_i)(t)\right|_{t=0}
=DF(\varphi(p))[(\varphi\circ\gamma_i)'(0)],
\qquad i=1,2.
\]
The coordinate velocities agree, so \(v_{[\gamma]}\)
is well defined. Linearity and the product rule show that
\[
v_{[\gamma]}(\lambda f+\mu g)
=\lambda v_{[\gamma]}(f)+\mu v_{[\gamma]}(g),
\qquad
v_{[\gamma]}(fg)
=f(p)v_{[\gamma]}(g)+g(p)v_{[\gamma]}(f).
\]
Therefore, \(v_{[\gamma]}\in T_pM\).

Write \(\varphi=(x^1,\dots,x^n)\). By
Lemma~\ref{lem:localidad-derivaciones-punto}, the coordinate functions
can be evaluated as germs, and
\[
(\varphi\circ\gamma)'(0)
=\bigl(v_{[\gamma]}(x^1),\dots,v_{[\gamma]}(x^n)\bigr).
\]
Thus \(v_{[\gamma_1]}=v_{[\gamma_2]}\) implies
\([\gamma_1]=[\gamma_2]\), and the map is injective.

For surjectivity, let \(v\in T_pM\) and set
\(x_0=\varphi(p)\) and
\[
w^i:=v(x^i),\qquad w=(w^1,\dots,w^n).
\]
Since \(p\) is interior, for sufficiently small \(|t|\) the curve
\[
\gamma(t):=\varphi^{-1}(x_0+tw),
\]
is defined, and \((\varphi\circ\gamma)'(0)=w\). Let \(f\in C^\infty(M)\) and write
\(F=f\circ\varphi^{-1}\).
Apply Taylor's theorem with integral remainder at the point $x_{0}$:
\[
R_{1,x_{0},F}(x)
= \displaystyle\sum_{|\alpha|=2}\frac{2}{\alpha !}(x-x_{0})^{\alpha}
\int_{0}^{1} (1-t) D^{\alpha}F\big(x_{0}+t(x-x_{0})\big)dt,
\]
and therefore
\[
F(x)
=
F(x_{0})
+\displaystyle\sum_{i=1}^{n}\frac{\partial F}{\partial x_{i}}(x_{0})(x_{i}-x_i^{(0)})
+R_{1,x_{0},F}(x).
\]
Taking $x=\varphi(q)$ for $q$ near $p$, we obtain:
\[
f(q)
=
f(p)
+\displaystyle\sum_{i=1}^{n}\frac{\partial F}{\partial x_{i}}(x_{0})\big(x^{i}(q)-x^{i}(p)\big)
+\underbrace{R_{1,x_{0},F}\big(\varphi(q)\big)}_{=:~r(q)}.
\]
We can rewrite this as
\[
r(q)
=
\displaystyle\sum_{|\alpha|=2}
\big(\varphi(q)-\varphi(p)\big)^{\alpha}h_{\alpha}(q),
\qquad
h_{\alpha}(q):=\frac{2}{\alpha !}\int_{0}^{1}(1-t)D^{\alpha}F\big(\varphi(p)+t(\varphi(q)-\varphi(p))\big)dt,
\]
where each $h_{\alpha}$ is smooth and well defined near $p$.

The vector $v$ annihilates the remainder term $r$.
For each multi-index $\alpha$ with $|\alpha|=2$, let
\[
m_{\alpha}(q):=\big(\varphi(q)-\varphi(p)\big)^{\alpha}=\prod_{i=1}^{n}\big(x^{i}(q)-x^{i}(p)\big)^{\alpha_{i}}.
\]
Then $m_{\alpha}(p)=0$, and $m_{\alpha}$ can be viewed as a product of at least two factors vanishing at $p$, since $|\alpha|=2$ forces either two of the $\alpha_{j}$ to be $1$ or one to be $2$.
By repeated application of the Leibniz rule,
\[
v\big(m_{\alpha}\big)=0
\quad\text{and}\quad
v\big(m_{\alpha}h_{\alpha}\big)=m_{\alpha}(p)v(h_{\alpha})+h_{\alpha}(p)v(m_{\alpha})=0.
\]
Summing over $|\alpha|=2$, we conclude that $v(r)=0$.
Applying $v$ to the decomposition of $f$ gives
\[
v(f)
=
\displaystyle\sum_{i=1}^{n}\frac{\partial F}{\partial x^{i}}(x_{0})v\big(x^{i}\big).
\]
On the other hand,
\[
v_{[\gamma]}(f)
=
\left.\frac{d}{dt}(f\circ\gamma)(t)\right|_{t=0}
=
\left.\frac{d}{dt}F(x_{0}+tw)\right|_{t=0}
=
\displaystyle\sum_{i=1}^{n}w^{i}\frac{\partial F}{\partial x^{i}}(x_{0}),
\]
and since $w^{i}=v(x^{i})$, we have $v(f)=v_{[\gamma]}(f)$ for every $f\in C^{\infty}(M)$.
Thus $v=v_{[\gamma]}$, and every derivation arises from an equivalence class of curves. We conclude that the map
\[
\{\text{classes of two-sided curves at $p$}\}
\longrightarrow T_pM,
\qquad [\gamma] \longmapsto v_{[\gamma]},
\]
is bijective when $p\in\operatorname{Int}(M)$.
\end{proof}

At a boundary point this equivalence no longer holds for two-sided
curves contained in $M$. We will return to this distinction after
describing the coordinate basis of $T_pM$.

There are other equivalent descriptions of the tangent space. One uses
\textit{germs of smooth functions}, which reflect the local nature of differentiation.

\begin{definition}\label{espacio de germenes}\index{germ of a smooth function@germ of a smooth function}
 Let $M$ be a smooth manifold with or without boundary and let $p\in M$.
 Define an equivalence relation on $C^{\infty}(M)$ by
 \[
 f\sim g
 \quad\Longleftrightarrow\quad
 \text{there exists a neighborhood $U$ of $p$ such that $f\restriction_{U}=g\restriction_{U}$}.
 \]
 We denote the equivalence class of $f$ by $[f]$ and call it
 the \textit{germ of $f$ at $p$}. The set of all such classes is
 denoted by $C_{p}^{\infty}(M)$ and is called the
 \textit{ring of germs of smooth functions at $p$}.
\end{definition}
\begin{proposition}[Derivations of the ring of germs]
\label{ejer:nociones-fundamentales-variedad-suave-frontera-recordemos-denota-anillo}
Let $M$ be a smooth manifold with or without boundary and let $p\in M$. There is a
canonical linear isomorphism
\[
 T_pM\longrightarrow
 \operatorname{Der}\bigl(C_p^\infty(M),\mathbb R\bigr),
 \qquad
 v\longmapsto\bigl([f]\longmapsto v(f)\bigr).
\]
Its inverse sends a derivation $D$ of germs to the derivation
$v_D(f):=D([f])$ of $C^\infty(M)$.
\end{proposition}
\begin{proof}
Addition and multiplication of germs are defined by
$[f]+[g]=[f+g]$ and $[f][g]=[fg]$; these operations are independent of
representatives because equality of germs is preserved by addition and
multiplication of functions. Moreover, $f(p)$ depends only on $[f]$. Therefore,
the rule
\[
 D([f][g])=f(p)D([g])+g(p)D([f])
\]
is well defined.

If $v\in T_pM$, Lemma~\ref{lem:localidad-derivaciones-punto} shows that
$\widetilde v([f]):=v(f)$ is independent of the representative. Linearity and the
Leibniz rule for $\widetilde v$ follow from those for $v$. Conversely, if
$D\in\operatorname{Der}(C_p^\infty(M),\mathbb R)$, then
$v_D(f):=D([f])$ is linear and
\[
 v_D(fg)=D([f][g])=f(p)v_D(g)+g(p)v_D(f),
\]
so $v_D\in T_pM$. Finally,
$v_{\widetilde v}(f)=\widetilde v([f])=v(f)$ and
$\widetilde{v_D}([f])=v_D(f)=D([f])$. The two constructions are inverse to each other
and require no choices.
\end{proof}

 With the tangent space defined, we can now introduce the tangent map, or differential, of a smooth function.

 \begin{definition}\label{def:nociones-fundamentales-mapeo-tangente-o-diferencial-de-en}\index{tangent map or differential} Let $M$ and $N$ be smooth manifolds with or without boundary, and let $F\colon M\longrightarrow N$ be a smooth map. For each $p\in M$, define the map \[dF_{p}\colon T_{p}M\longrightarrow T_{F(p)}N\] called the \textbf{tangent map or differential of $F$ at $p$}, by
 \[
 dF_p(v)(f):=v(f\circ F),
 \qquad v\in T_pM,\quad f\in C^\infty(N).
 \]
Indeed, $dF_p(v)$ is linear and, for $f,g\in C^\infty(N)$,
\[
 dF_p(v)(fg)
 =f(F(p))dF_p(v)(g)+g(F(p))dF_p(v)(f),
\]
so it belongs to $T_{F(p)}N$.
 \end{definition}
 The tangent map, or differential, has the following properties:
 \begin{proposition}
 \label{propiedadesmapeotangente}
 Let $M$, $N$, and $P$ be smooth manifolds with or without boundary, and let $F\colon M\longrightarrow N$ and $G\colon N\longrightarrow P$ be smooth maps. Let $p\in M$.
 \begin{enumerate}[label=(\alph*)]
 \item $dF_{p}\colon T_{p}M\longrightarrow T_{F(p)}N$ is linear.
 \item $d(G\circ F)_{p}=dG_{F(p)}\circ dF_{p}\colon T_{p}M\longrightarrow T_{G(F(p))}P$.
 \item $d(\operatorname{Id}_{M})_{p}=\operatorname{Id}_{T_{p}M}$.
 \item If $F$ is a diffeomorphism, then $dF_{p}\colon T_{p}M\longrightarrow T_{F(p)}N$ is an isomorphism and \[(dF_{p})^{-1}=d(F^{-1})_{F(p)}.\]
 \end{enumerate}
 \end{proposition}
\begin{proof}
Linearity follows directly from the linearity of each derivation.
If $v\in T_pM$ and $h\in C^\infty(P)$, then
\[
 d(G\circ F)_p(v)(h)
 =v(h\circ G\circ F)
 =dF_p(v)(h\circ G)
 =dG_{F(p)}(dF_p(v))(h),
\]
which proves the chain rule. The formula for the identity is the
immediate case $h\circ\operatorname{Id}_M=h$. If $F$ is a
diffeomorphism, apply the chain rule to
$F^{-1}\circ F=\operatorname{Id}_M$ and to
$F\circ F^{-1}=\operatorname{Id}_N$; thus
$d(F^{-1})_{F(p)}$ is the inverse of $dF_p$.
\end{proof}
 Since $T_{p}M$ is a vector space, a basis is useful for describing it. In the Euclidean case, $T_{p}\mathbb{R}^{n}$ is isomorphic to $\mathbb{R}^{n}$ for every $p\in\mathbb{R}^{n}$, and tangent vectors may be viewed as directional derivatives at $p$ in the direction of $v\in T_{p}\mathbb{R}^{n}$. This leads to the following proposition:
 \begin{proposition}

 \label{basedetpm} Let $M$ be a smooth manifold with or without boundary of dimension $n$, and let $p\in M$. Then $T_{p}M$ is a vector space of dimension $n$, and for each smooth chart $(U,\phi)$ such that $p\in U$, with $\phi=(x^{1},\dots,x^{n})$, the vectors
 \[
 \left.\boldsymbol{\partial}_{i}\right|_{p}(f)
 :=\frac{\partial \widetilde F}{\partial x^{i}}(\phi(p)),
 \qquad F=f\circ\phi^{-1},
 \]
 form a basis for $T_pM$. Here $\widetilde F=F$ in an interior chart;
 in a boundary chart, $\widetilde F$ is any smooth extension of
 $F$ to a Euclidean neighborhood of $\phi(p)$.

 \end{proposition}

\begin{proof}
In a boundary chart, the definition of a smooth function provides the
extension $\widetilde F$. If two extensions agree on the relatively open set
$\phi(U)\subseteq\mathbb H^n$, their partial derivatives agree
at interior points and, by continuity, also at boundary points. Thus
the formula is independent of the extension. Euclidean linearity and the
Leibniz rule show that each
$\left.\boldsymbol{\partial}_{i}\right|_p$ is a derivation at $p$.

Let $v\in T_pM$ and write $x_0=\phi(p)$. Restrict the chart, if
necessary, so that segments starting at $x_0$ remain in
its image. Hadamard's lemma~\ref{lem:hadamard-finito-dimensional},
applied to a Euclidean extension when $p\in\partial M$, gives smooth
functions $F_i$ such that, near $x_0$,
\[
 F(x)-F(x_0)=\displaystyle\sum_{i=1}^n(x^i-x_0^i)F_i(x),
 \qquad F_i(x_0)=\frac{\partial\widetilde F}{\partial x^i}(x_0).
\]
Lemma~\ref{lem:localidad-derivaciones-punto} allows us to apply $v$ to this
identity of germs. Since each $x^i-x_0^i$ vanishes at $p$, we obtain
\[
 v(f)=\displaystyle\sum_{i=1}^n v(x^i)
 \frac{\partial\widetilde F}{\partial x^i}(x_0)
 =\left(\displaystyle\sum_{i=1}^n v(x^i)
 \left.\boldsymbol{\partial}_{i}\right|_p\right)(f).
\]
Thus the coordinate vectors span $T_pM$. If
$\displaystyle \sum_{i=1}^{n} a^i\left.\boldsymbol{\partial}_{i}\right|_p=0$, evaluation on the
germ of $x^j$ gives $a^j=0$. They are therefore linearly independent and
form a basis; in particular, $\dim T_pM=n$.

If $(y^1,\ldots,y^n)$ is another coordinate system around $p$, the
expansion just proved, applied to
$v=\left.\boldsymbol{\partial}_{i}\right|_p$, gives the change-of-basis
formula
\begin{equation}
\label{eq:cambio-base-vectores-coordenados}
 \left.\boldsymbol{\partial}_{i}\right|_p
 =\displaystyle\sum_{j=1}^n
 \frac{\partial y^j}{\partial x^i}(p)
 \left.\boldsymbol{\partial}_{j}\right|_p.
\end{equation}
\end{proof}

\begin{remark}[Two-sided curves at the boundary]
\label{obs:curvas-bilaterales-frontera-tangentes}
Suppose that $p\in\partial M$ and take a boundary chart with
$\mathbb H^n=\{x^n\geq0\}$ and $x^n(p)=0$. If
$\gamma\colon(-\varepsilon,\varepsilon)\to M$ is smooth and
$\gamma(0)=p$, the function $t\mapsto x^n(\gamma(t))$ is nonnegative and
attains a minimum at $0$; therefore,
$(x^n\circ\gamma)'(0)=0$. The velocities of two-sided curves contained
in $M$ range over precisely the coordinate subspace
\[
 \operatorname{span}\left\{
 \left.\boldsymbol{\partial}_{1}\right|_p,\ldots,
 \left.\boldsymbol{\partial}_{n-1}\right|_p
 \right\},
\]
which will later be identified with $T_p\partial M$, rather than all of $T_pM$.
The reverse inclusion follows by using the coordinate
curves $t\mapsto x_0+t(w^1,\ldots,w^{n-1},0)$.

With this convention, a vector
$v=v^i\left.\boldsymbol{\partial}_{i}\right|_p$ points inward if
$v^n>0$, outward if $v^n<0$, and is tangent to the boundary if $v^n=0$.
The sign is independent of the chart: for two normal boundary coordinates
$x^n$ and $y^n$, Hadamard's lemma gives
$y^n=a x^n$ in a relative neighborhood of $p$, with $a(p)>0$. For this reason, the primary definition
of $T_pM$ at a boundary point is the definition by
derivations, rather than the description by two-sided curves.
\end{remark}

\begin{proposition}[Tangent spaces and velocities at corners]
\label{prop:esquinas-tangente-conos}
Let $M$ be a smooth manifold with corners of dimension $n$, and let
$p\in M$ have depth $d$. The definition by derivations gives a space
$T_pM$ of dimension $n$; the coordinate vectors
$\left.\boldsymbol{\partial}_i\right|_p$ form a basis.
In a chart reduced to $G_{n,d}$ at $p$, velocities of one-sided
curves form the cone $C_p$ of nonnegative normal coordinates, whereas
velocities of two-sided curves form $C_p\cap(-C_p)$, of dimension $n-d$.
\end{proposition}

\begin{proof}
Lemma~\ref{lem:esquinas-jets-extension-local} makes partial derivatives
of a smooth function at $p$ independent of the extension. They define
derivations by the Leibniz rule. Locality of derivations also holds
here: if $f=0$ near $p$, choose a cutoff $\chi$ supported there and
equal to one near $p$, by
Proposition~\ref{prop:esquinas-particiones-unidad}. Then $f\chi=0$,
and the Leibniz rule gives $v(f)=v(f\chi)=0$. This allows evaluation
on germs of local functions using cutoffs.

Choose a Euclidean extension $\widetilde F$ of
$F=f\circ\phi^{-1}$ on a convex ball centered at $z=\phi(p)$.
The fundamental theorem of calculus gives
\[
 F(x)-F(z)=\sum_{i=1}^n(x^i-z^i)F_i(x),\qquad
 F_i(x):=\int_0^1
       \partial_i\widetilde F(z+t(x-z))\,dt.
\]
These functions are smooth by differentiation under the integral, and
$F_i(z)=\partial_i\widetilde F(z)$. Applying a derivation $v$ to the
identity of germs gives
\[
 v(f)=\sum_{i=1}^n v(x^i)\partial_i\widetilde F(z).
\]
Thus the coordinate vectors span $T_pM$. Evaluation on the germs of
$x^j$ gives the identity matrix, proving linear independence.
The chain rule gives the same change of basis as in
\eqref{eq:cambio-base-vectores-coordenados}; the Jacobians are
invertible by invariance of depth. Hence these bases also define
the ordinary smooth tangent bundle.

The description of the one-sided cone was proved in
Proposition~\ref{prop:esquinas-profundidad-invariante}. For a two-sided
curve, each active normal coordinate has a minimum at zero, so its
derivative vanishes. Conversely, a coordinate line whose active
normal components vanish stays in the model for sufficiently small
times of both signs. Its velocity realizes any vector in
$C_p\cap(-C_p)$. Exactly $d$ coordinates vanish, giving dimension $n-d$.
\end{proof}

Vector fields and vector bundles use this full tangent space, rather
than only two-sided velocities. The components of a field or section
are continuous, of class $C^k$, or smooth if and only if they have
that regularity in every local frame in the corner model. Indeed,
coordinates and trivializations are homeomorphisms, and component
changes multiply by smooth invertible matrices; the Leibniz rule and
composition with chart transitions preserve each of these classes.
This permits bundle metrics and connections smooth up to all faces
with the same local expressions as in the case without corners.

 This notation is suggestive because it essentially transfers the concept of a partial derivative from the Euclidean setting to manifolds. The following definition makes this precise:
\begin{definition}\label{def:nociones-fundamentales-derivada-carta-suave-funcion-suave}\index{partial derivative!on a manifold}
Let $M$ be a smooth manifold with or without boundary, let $(U,\phi)$ be a smooth
chart of $M$, with $\phi=(x^{1},\dots,x^{n})$, and let
$f\colon U\longrightarrow\mathbb R$ be a smooth function. Given a multi-index
$\alpha=(\alpha_{1},\dots,\alpha_{n})\in\mathbb N_0^{n}$, we define:

\begin{enumerate}[label=(\alph*)]
 \item For each $p\in U$, the partial derivative of $f$ associated with $\alpha$ at $p$ by
 $\partial^{\alpha}f(p):=\partial^{\alpha}(f\circ\phi^{-1})(\phi(p))$.

 \item The partial derivative function of $f$ associated with $\alpha$ on $U$ by
 $\partial^{\alpha}f\colon U\longrightarrow\mathbb R$, given by $p\mapsto \partial^{\alpha}f(p)$.
\end{enumerate}
\end{definition}
 The tangent bundle, which gives this section its name, is defined as follows:
 \begin{definition}\label{def:nociones-fundamentales-haz-tangente-de}\index{tangent bundle}
 Let $M$ be a smooth manifold with or without boundary. Define the \textbf{tangent bundle of $M$} by \[TM:=\coprod_{p\in M}T_{p}M\]
 \end{definition}
\noindent\textbf{Index convention.}
In coordinates $(x^1,\ldots,x^n)$, coordinate indices range over
$\{1,\ldots,n\}$. Under the Einstein convention, an index appearing
exactly twice in one term, once as an upper index and once as a lower
index, indicates summation over its full range. For example,
\[
v^i\boldsymbol\partial_i
=\displaystyle\sum_{i=1}^n v^i\boldsymbol\partial_i,
\qquad
\omega_i v^i=\displaystyle\sum_{i=1}^n\omega_i v^i.
\]
Free indices remain fixed, and a component identity holds for each
choice of them. A fiber index of a rank-$r$ bundle ranges over
$\{1,\ldots,r\}$; when different spaces occur, their respective
index ranges will be specified. Explicit sums and Einstein notation
represent the same operation.

 Elements of $TM$ have the form $(p,v)$ with $p\in M$, $v\in T_{p}M$. The following notations are commonly used for an element of $TM$: $(p,v)$, $v_{p}$, and $v$, the last when we do not wish to emphasize the base point. By Proposition~\ref{basedetpm}, every element of $T_{p}M$ can be written in the form (using Einstein notation) \[v_{p}=v^{i}\boldsymbol{\partial}_{i}\biggr|_{p},\] which gives another way to express elements of $TM$.

 The tangent bundle can be equipped with a smooth manifold structure, with or without boundary according to the type of manifold we start with. The following proposition describes the construction:
 \begin{proposition}\label{estructurasuavedelhaztangente}
 If $M$ is a smooth $n$-manifold with or without boundary, the tangent bundle $TM$ has a natural topology and smooth structure making it a smooth $2n$-manifold with or without boundary. With respect to this smooth structure, the projection $\pi\colon TM\longrightarrow M$ is smooth. The coordinate charts for this smooth structure are given as follows:
 \begin{enumerate}[label=(\alph*)]
 \item If $(U,\phi)$ is an interior chart of $M$, then $\widetilde{\phi}\colon \pi^{-1}(U)\longrightarrow \mathbb{R}^{2n}$, given by \[\widetilde{\phi}\left(p,v^{i}\boldsymbol{\partial}_{i}\biggr|_{p}\right)=(x^{1}(p),\dots,x^{n}(p),v^{1},\dots,v^{n})\] makes $(\pi^{-1}(U),\widetilde{\phi})$ an interior chart of $TM$.
 \item If $(U,\phi)$ is a boundary chart of $M$, then $\widetilde{\phi}\colon \pi^{-1}(U)\longrightarrow \mathbb{H}^{2n}$, given by \[\widetilde{\phi}\left(p,v^{i}\boldsymbol{\partial}_{i}\biggr|_{p}\right)=(v^{1},\dots,v^{n},x^{1}(p),\dots,x^{n}(p))\] makes $(\pi^{-1}(U),\widetilde{\phi})$ a boundary chart of $TM$.
 \end{enumerate}\end{proposition}
\begin{proof}
For each chart $(U,\phi)$, the map $\widetilde\phi$ is a bijection
onto $\phi(U)\times\mathbb R^n$ in the interior case and onto
$\mathbb R^n\times\phi(U)$ in the boundary case; the latter set is
open in $\mathbb H^{2n}$ under the convention in the statement. If
$(V,\psi)$ is another chart, write $y=\psi\circ\phi^{-1}(x)$. The change-of-basis
formula \eqref{eq:cambio-base-vectores-coordenados} shows that,
on the overlap, the components of a vector satisfy
\[
 w^j=\displaystyle\sum_{i=1}^n\frac{\partial y^j}{\partial x^i}(x)v^i.
\]
Thus, apart from explicitly interchanging the order of the two blocks in
a boundary chart, the change of coordinates on $TM$ is
\[
 (x,v)\longmapsto\bigl(y(x),Dy(x)v\bigr).
\]
It is smooth, and its inverse is obtained by replacing $y$ by $x(y)$. An
overlap between an interior chart and a boundary chart contains no points
of $\partial M$, so the same formula is interpreted on Euclidean open
sets.

The smooth chart lemma proved above, applied to this family,
provides a unique Hausdorff, second countable topology and the
stated smooth structure on $TM$. In these coordinates, the projection is
$(x,v)\mapsto x$ (with the corresponding interchange of blocks at the
boundary), so it is smooth. If $M$ has no boundary, all charts are
interior charts; if it has boundary, the boundary charts described above determine
the boundary of $TM$. This also completes the uniqueness of the structure
for which the charts in the statement are smooth.
\end{proof}
\section{Immersions, submersions, and embeddings}

As in the Euclidean case, the differential of a function between manifolds, which plays the role of the usual derivative, describes the best linear approximation to a smooth map near each point. Thus many local properties of a map can be studied through the linear properties of its differential. This motivates the following definition.

\begin{definition}\label{def:nociones-fundamentales-variedad-5}\index{rank of a smooth map}
Let $M$ and $N$ be smooth manifolds with or without boundary, and let $F\colon M\longrightarrow N$ be a smooth map. We define the rank of $F$ at a point $p\in M$ as the rank of the linear map
\[
dF_{p}\colon T_{p}M\longrightarrow T_{F(p)}N.
\]
Equivalently, the rank of $F$ at $p$ is the dimension of the image of $dF_{p}$. In coordinates, this number equals the rank of the Jacobian matrix of a coordinate representation of $F$. If $F$ has the same rank at every point of its domain, we say that $F$ has constant rank.
\end{definition}

The notions of immersion and submersion generalize, respectively, injectivity and surjectivity of linear maps to the smooth setting.

\begin{definition}\label{def:inmersion-submersion}\index{immersion}\index{submersion}
Let $M$ and $N$ be smooth manifolds with or without boundary, and let
$F\colon M\longrightarrow N$ be a smooth map. We say that:
\begin{enumerate}[label=(\alph*)]
 \item $F$ is an immersion if $dF_{p}$ is injective for every $p\in M$;

 \item $F$ is a submersion if $dF_{p}$ is surjective for every $p\in M$.
\end{enumerate}
\end{definition}

Being an immersion or a submersion is a local property. If the differential of a map is injective or surjective at a point, the same property holds throughout some open neighborhood of that point.

\begin{proposition}\label{prop:nociones-fundamentales-mapeo-suave-label-inyectivo-vecindad-abierta}
Let $M$ and $N$ be smooth manifolds without boundary, let
$F\colon M\longrightarrow N$ be a smooth map, and let $p\in M$.

\begin{enumerate}[label=(\alph*)]
 \item If $dF_{p}$ is injective, there exists an open neighborhood $U$ of $p$ such that $F\restriction_{U}$ is an immersion.

 \item If $dF_{p}$ is surjective, there exists an open neighborhood $U$ of $p$ such that $F\restriction_{U}$ is a submersion.
\end{enumerate}
\end{proposition}
\begin{proof}
In charts around \(p\) and \(F(p)\), the differential is represented
by the Jacobian matrix. If \(dF_p\) is injective and \(m=\dim M\), some
minor of order \(m\) has nonzero determinant. By continuity,
this minor remains nonsingular in a neighborhood of \(p\); there the
rank is \(m\), so the differential is injective. If \(dF_p\) is
surjective and \(n=\dim N\), the same argument applies to a minor of
order \(n\), and the differential is surjective in a neighborhood of \(p\).
\end{proof}

We now give some examples of immersions and submersions.

\begin{example}\label{ej:nociones-fundamentales-label-variedades-suaves-proyeccion-submersion-efecto}
\begin{enumerate}[label=(\alph*)]

\item If $M_{1},\dots,M_{k}$ are smooth manifolds without boundary, each projection
$\pi_{i}\colon M_{1}\times\cdots\times M_{k}\longrightarrow M_{i}$ is a submersion. For each $(p_{1},\dots,p_{k})\in M_{1}\times\cdots\times M_{k}$, its differential satisfies
$d(\pi_{i})_{(p_{1},\dots,p_{k})}(v_{1},\dots,v_{k})=v_{i}$. Indeed, given
$w\in T_{p_i}M_i$, the vector $(0,\dots,0,w,0,\dots,0)$ in the tangent space
to the product is sent to $w$; therefore, the differential is
surjective.

\item The map $F\colon \mathbb R^{m}\longrightarrow\mathbb R^{n}$, given by
$F(x^{1},\dots,x^{m})=(x^{1},\dots,x^{m},0,\dots,0)$, with $m\leq n$, is an immersion. Its differential is the natural linear inclusion $\mathbb R^{m}\hookrightarrow\mathbb R^{n}$, which is injective.

\item The map $G\colon \mathbb R^{3}\longrightarrow\mathbb R^{2}$, given by $G(x,y,z)=(x,y)$, is a submersion, since $dG_{(x,y,z)}(a,b,c)=(a,b)$ for every $(x,y,z)\in\mathbb R^{3}$, and this linear map is surjective.

\item The map $\gamma\colon \mathbb R\longrightarrow\mathbb R^{2}$, given by $\gamma(t)=(\cos t,\sin t)$, is an immersion, since $\gamma'(t)=(-\sin t,\cos t)\neq 0$ for every $t\in\mathbb R$.

\item The map $H\colon \mathbb R^{2}\longrightarrow\mathbb R$, given by $H(x,y)=x^{2}+y^{2}$, is not a submersion on all of $\mathbb R^{2}$, since $dH_{(0,0)}=0$. It is, however, a submersion on $\mathbb R^{2}\setminus\{(0,0)\}$.

\item Consider the map $K\colon \mathbb R^{2}\longrightarrow\mathbb R^{2}$, given by $K(x,y)=(x^{2}-y^{2},2xy)$. Identifying $\mathbb R^{2}$ with $\mathbb C$, this map corresponds to $z\mapsto z^{2}$. A direct calculation shows that
$DK(x,y)=\begin{pmatrix}2x&-2y\\2y&2x\end{pmatrix}$, whose determinant is $4(x^{2}+y^{2})$. Thus $K$ is both an immersion and a submersion at every point other than the origin, but has neither property at $(0,0)$.

\end{enumerate}
\end{example}
The most important theorems of multivariable differential calculus also hold for smooth manifolds and will be useful throughout the book. We present them next.

\begin{theorem}[inverse function theorem]\label{teo: funcion inversa variedades}\index{inverse function theorem@inverse function theorem}
Let $M$ and $N$ be smooth manifolds without boundary, and let $F\colon M\longrightarrow N$ be a smooth map. If $p\in M$ is such that $dF_p$ is invertible, there exist connected open neighborhoods $U_0$ of $p$ and $V_0$ of $F(p)$ such that $F\restriction_{U_0}\colon U_0\longrightarrow V_0$ is a diffeomorphism.
\end{theorem}

\begin{proof}
Since $dF_p$ is a bijection, we have $\dim(M)=\dim(N)$. Consider smooth charts $(U,\phi)$ and $(V,\psi)$ centered at $p$ and $F(p)$, respectively, with $F(U)\subseteq V$. Then $\widetilde F:=\psi\circ F\circ\phi^{-1}$ maps the open set $\widetilde U:=\phi(U)\subseteq\mathbb R^n$ into the open set $\widetilde V:=\psi(V)\subseteq\mathbb R^n$.

Since $\phi$ and $\psi$ are diffeomorphisms, the differential
\[
d\widetilde F_{\phi(p)}
=
d\psi_{F(p)}\circ dF_p\circ d(\phi^{-1})_{\phi(p)}
\]
is an isomorphism.
Theorem~\ref{lem:funcion-inversa-euclidea-previa-variedades}, applied to
$\widetilde F$, gives connected open sets
$\widetilde U_0\subseteq\widetilde U$ and
$\widetilde V_0\subseteq\widetilde V$, with
$\phi(p)\in\widetilde U_0$ and
$\widetilde F(\phi(p))\in\widetilde V_0$, such that
$\widetilde F\restriction_{\widetilde U_0}\colon
\widetilde U_0\longrightarrow\widetilde V_0$ is a diffeomorphism.

Defining $U_0:=\phi^{-1}(\widetilde U_0)$ and $V_0:=\psi^{-1}(\widetilde V_0)$, we obtain connected open neighborhoods of $p$ and $F(p)$, respectively, such that $F\restriction_{U_0}\colon U_0\longrightarrow V_0$ is a diffeomorphism.
\end{proof}

\begin{remark}\label{obs:nociones-fundamentales-enunciado-anterior-requiere-variedades-frontera-efecto}
The preceding statement requires the manifolds to have no boundary. If boundary is allowed, the result may fail. For example, the inclusion $\iota\colon \mathbb H^n\hookrightarrow\mathbb R^n$ has invertible differential at every point but is not a local diffeomorphism at points of $\partial\mathbb H^n$.
\end{remark}

This motivates the concept of a local diffeomorphism:
\begin{definition}\label{def:nociones-fundamentales-variedad-6}\index{local diffeomorphism}
 Let $M$ and $N$ be manifolds with or without boundary. A function $F\colon M\longrightarrow N$ is called a local diffeomorphism if every $p\in M$ has an open neighborhood $U$ such that $F(U)$ is open in $N$ and $F\restriction_{U}\colon U\longrightarrow F(U)$ is a diffeomorphism.
\end{definition}
Local diffeomorphisms have the following useful property.
\begin{proposition}\label{prop:nociones-fundamentales-variedades-suaves-frontera-funcion-difeomorfismo-local}
 Let $M$ and $N$ be smooth manifolds without boundary, and let $F\colon M\longrightarrow N$ be a function.
 \begin{enumerate}
 \item $F$ is a local diffeomorphism if and only if it is both a smooth immersion and a smooth submersion.
 \item If $\operatorname{dim}(M)=\operatorname{dim}(N)$ and $F$ is a smooth immersion or a smooth submersion, then it is a local diffeomorphism.
 \end{enumerate}
\end{proposition}
\begin{proof}
If \(F\) is a local diffeomorphism, it has a smooth inverse
around each \(p\in M\). The chain rule shows that \(dF_p\) and
\(d(F^{-1})_{F(p)}\) are inverse to each other; thus \(F\) is an immersion and
a submersion.

Conversely, if \(F\) is an immersion and a submersion, each \(dF_p\) is an
isomorphism. Theorem~\ref{teo: funcion inversa variedades} provides
neighborhoods on which \(F\) is a diffeomorphism, so \(F\) is a
local diffeomorphism. This proves (1). When the dimensions are equal, every
injective linear map between the tangent spaces is
surjective, and conversely; (2) follows from (1).
\end{proof}
The rank theorem (Theorem~\ref{teo: del rango euclidiano}) has an analogue for smooth manifolds, bearing the same name.
\begin{theorem}[Rank theorem]\label{teo: teo del rango}\index{rank theorem}
Let $M$ and $N$ be smooth manifolds without boundary of dimensions $m$ and $n$, respectively, and let $F\colon M\longrightarrow N$ be a smooth map of constant rank $r$. Then, for every $p\in M$, there exist smooth charts $(U,\phi)$ of $M$ centered at $p$ and $(V,\psi)$ of $N$ centered at $F(p)$ such that $F(U)\subseteq V$ and the coordinate representation $\widetilde F:=\psi\circ F\circ\phi^{-1}$ has the form $\displaystyle \widetilde F(x^1,\dots,x^m)=(x^1,\dots,x^r,0,\dots,0)$.

In particular, if $F$ is an immersion, then $r=m$ and $\displaystyle \widetilde F(x^1,\dots,x^m)=(x^1,\dots,x^m,0,\dots,0)$. Likewise, if $F$ is a submersion, then $r=n$ and $\displaystyle \widetilde F(x^1,\dots,x^m)=(x^1,\dots,x^n)$.
\end{theorem}

\begin{proof}
Fix $p\in M$. Take smooth charts $(U_0,\phi_0)$ of $M$ around $p$ and $(V_0,\psi_0)$ of $N$ around $F(p)$, with $F(U_0)\subseteq V_0$. Consider the coordinate representation $\displaystyle f:=\psi_0\circ F\circ\phi_0^{-1}\colon \phi_0(U_0)\subseteq\mathbb R^m\longrightarrow \mathbb R^n$. Since $F$ has constant rank $r$, the derivative of $f$ also has constant rank $r$. Moreover, since $F$ is smooth, $f$ is of class $C^\infty$. Applying the smooth version of Theorem~\ref{teo: del rango euclidiano} gives open sets $\widetilde W_1,\widetilde W_2\subseteq\mathbb R^m$ and $W_1,W_2\subseteq\mathbb R^n$, together with smooth diffeomorphisms $\displaystyle \Phi\colon \widetilde W_2\longrightarrow \widetilde W_1$ and $\displaystyle \Psi\colon W_2\longrightarrow W_1$, such that
\[
\Psi^{-1}\circ f\circ \Phi(x^1,\dots,x^m)
=
(x^1,\dots,x^r,0,\dots,0).
\]

Define $\displaystyle U:=\phi_0^{-1}(\widetilde W_1)$ and $\displaystyle V:=\psi_0^{-1}(W_1)$. Restricting $U$ if necessary, we may assume that $F(U)\subseteq V$. Now take the charts $\displaystyle \phi:=\Phi^{-1}\circ \phi_0\restriction_U$ and $\displaystyle \psi:=\Psi^{-1}\circ \psi_0\restriction_V$. Then $(U,\phi)$ is a smooth chart of $M$ centered at $p$, and $(V,\psi)$ is a smooth chart of $N$ centered at $F(p)$. In these charts,
\[
\psi\circ F\circ\phi^{-1}
=
\Psi^{-1}\circ\psi_0\circ F\circ\phi_0^{-1}\circ\Phi
=
\Psi^{-1}\circ f\circ\Phi.
\]
Therefore, $\displaystyle \widetilde F(x^1,\dots,x^m)=(x^1,\dots,x^r,0,\dots,0)$.

If $F$ is an immersion, then $r=m$, and the preceding expression becomes $\displaystyle \widetilde F(x^1,\dots,x^m)=(x^1,\dots,x^m,0,\dots,0)$. If $F$ is a submersion, then $r=n$, and the expression reduces to $\displaystyle \widetilde F(x^1,\dots,x^m)=(x^1,\dots,x^n)$.
\end{proof}
\section{Submanifolds}

A submanifold allows us to treat a subset as a smooth space in its own right. The definition must reconcile two features: the induced topology and injectivity of the differential. The following figures motivate this condition before we formulate it precisely.
\begin{figura}[h]    \includegraphics[width=0.8\linewidth]{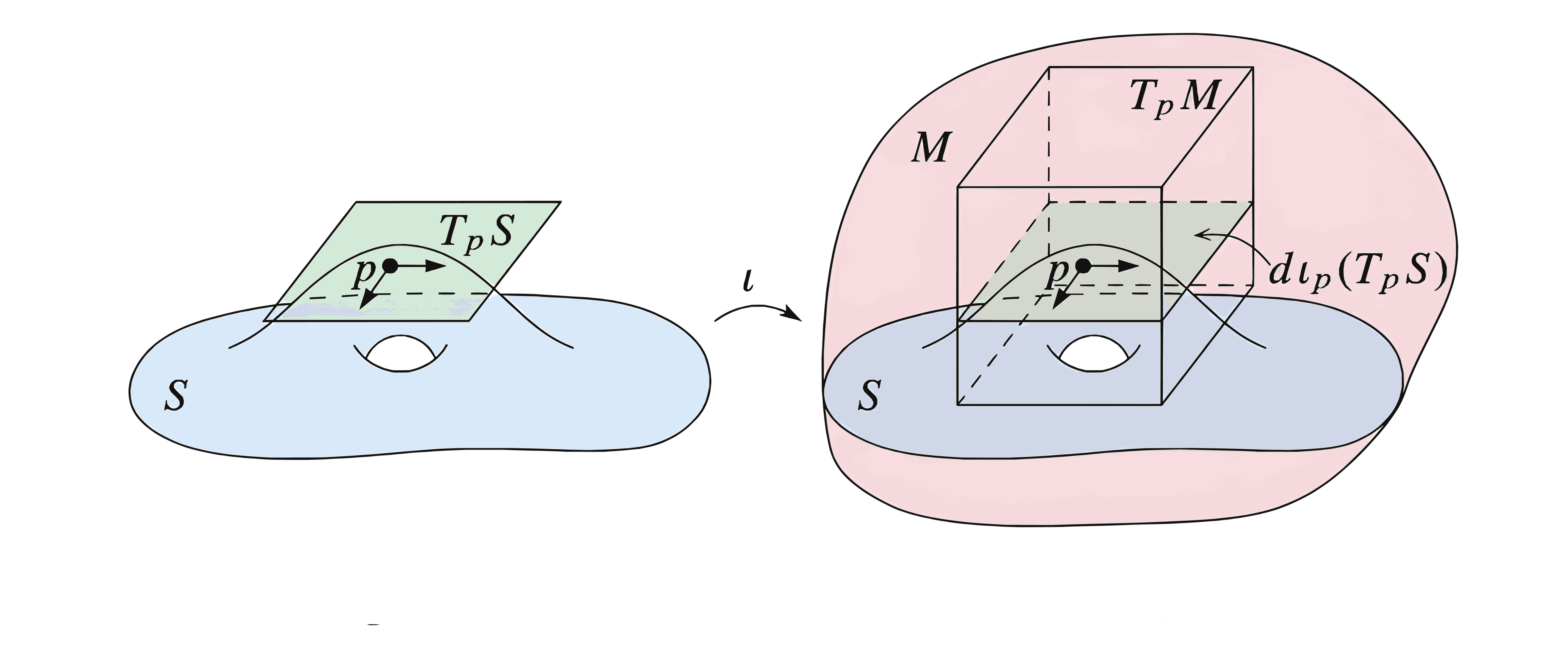}
    \caption{Tangent space to a submanifold.}
    \label{fig:tangente-subvariedad}
\end{figura}

\begin{definition}[Embeddings and submanifolds]
\label{def:subvariedades-inmersas-encajadas-finitas}
\index{smooth embedding}\index{immersed submanifold}\index{embedded submanifold}
Let $M$ and $N$ be smooth manifolds with or without boundary.
The manifold structures considered on subsets may also have boundary.
\begin{enumerate}[label=(\alph*)]
\item A smooth map $F\colon N\to M$ is a \textbf{smooth embedding} if it is
an injective immersion and a homeomorphism from $N$ onto $F(N)$ with the
topology induced by $M$.
\item A subset $S\subseteq M$, equipped with a smooth manifold
structure, is an \textbf{immersed submanifold} if the inclusion
$\iota\colon S\hookrightarrow M$ is an immersion.
\item It is an \textbf{embedded submanifold} if $\iota$ is a smooth embedding.
It is \textbf{properly embedded} if, in addition, $\iota$ is proper.
\end{enumerate}
For an immersed submanifold, the topology on $S$ may be finer than
the subspace topology; for an embedded submanifold, the two coincide.
\end{definition}

The Whitney embedding theorem ensures that every finite-dimensional
smooth manifold, under the Hausdorff and second-countability assumptions
adopted here, admits a smooth embedding into some Euclidean space;
see \cite[Theorem 6.15]{LeeS}. We only mention this result.
Existence of an embedding singles out no particular embedding and does
not provide a global chart. Functions, fields, and sections are defined
on the manifold and checked locally, independently of this choice.

\begin{definition}[Slice chart]
\label{def:carta-rebanada-subvariedad}
\index{slice chart}\index{chart adapted to a submanifold}
Let $M$ be a smooth manifold without boundary, and let $S\subseteq M$ and $p\in S$.
A chart $(U,\phi)$ of $M$ is
\textbf{adapted to $S$ at $p$}, or is a \textbf{slice chart}, if
$p\in U$, $\phi(p)=0$, and, for some $k$,
\[
 \phi(U\cap S)
 =\phi(U)\cap\bigl(\mathbb R^k\times\{0\}^{n-k}\bigr).
\]
The last equality includes the local topology of the subset, not merely
an equality of tangent spaces.
\end{definition}

\begin{theorem}[Local slice criterion]
\label{teo:criterio-cartas-rebanada-subvariedades}
\index{local slice criterion}
A subset $S$ of a smooth $n$-manifold without boundary admits the
structure of an embedded $k$-dimensional submanifold without boundary
if and only if every $p\in S$ has a $k$-dimensional slice chart.
When it exists, this smooth structure is unique.
\end{theorem}
\begin{proof}
Suppose such charts exist. Give $S$ the subspace topology and, for
each adapted chart $(U,\phi)$, define
\[
\phi_S:=\operatorname{pr}_{\mathbb R^k}\circ
\phi\restriction_{U\cap S}.
\]
The slice equality shows that $\phi_S$ is a homeomorphism onto an
open subset of $\mathbb R^k$. If $(V,\psi)$ is another adapted chart,
the transition map, on its domain, is
\[
u\longmapsto\operatorname{pr}_{\mathbb R^k}
\bigl((\psi\circ\phi^{-1})(u,0)\bigr).
\]
It is smooth, as is its inverse. Moreover, the subspace $S$ is
Hausdorff and second countable because $M$ has both properties.
These charts therefore define a smooth structure without boundary
on $S$. In these coordinates the inclusion is $u\mapsto(u,0)$,
whose differential is injective; the chosen topology makes it an
embedding.

Conversely, suppose $\iota\colon S\hookrightarrow M$ is a smooth
embedding. Theorem~\ref{teo: teo del rango}, applied to $\iota$
near $p$, gives charts $\theta\colon W_S\to A\subseteq\mathbb R^k$
and $\phi\colon U_0\to B\subseteq\mathbb R^n$, centered at $p$,
in which $\phi(\iota(\theta^{-1}(u)))=(u,0)$ for $u\in A$;
shrink $W_S$ so that its image lies in $U_0$. Since $S$ has the
subspace topology, there is an open subset $O$ of $M$ such that
$p\in O$ and $O\cap S\subseteq W_S$. Choose open boxes
$A_1\subseteq A$ and $C_1\subseteq\mathbb R^{n-k}$, both centered
at zero, with
\[
A_1\times C_1\subseteq\phi(O\cap U_0).
\]
The set $U:=\phi^{-1}(A_1\times C_1)$ contains no points of $S$
outside $W_S$. Conversely, every $(u,0)$ with $u\in A_1$
corresponds to $\theta^{-1}(u)\in S$. Thus
$\phi(U\cap S)=A_1\times\{0\}$, as required. The same construction
shows that the restricted charts are smooth for any structure
making $\iota$ an embedding. Since they cover $S$, they determine
its maximal atlas and prove uniqueness.
\end{proof}

\begin{proposition}[Maps into an embedded submanifold]
\label{prop:factorizacion-suave-subvariedad}
Let $S$ be an embedded submanifold without boundary of a smooth
manifold $M$ without boundary. If $F\colon N\to M$ is smooth,
where $N$ may have boundary, and $F(N)\subseteq S$, then the same
map with codomain $S$ is smooth.
\end{proposition}
\begin{proof}
The map into $S$ is continuous by the subspace topology. In a slice
chart $\phi$ around $F(p)$, its coordinates in $S$ are the first
$\dim S$ components of $\phi\circ F$. They are smooth near $p$;
if $p\in\partial N$, they admit the smooth extensions required by
the definition of smoothness on a half-space. This proves the claim
at every point of $N$.
\end{proof}

The embedding condition matters: it recovers smoothness into $S$
from ambient smoothness. For a merely immersed submanifold, even
continuity into its finer topology need not follow from continuity
into $M$.

\begin{proposition}[Local embeddings and proper embeddings]
\label{prop:inmersion-local-encaje-propio}
Let $N$ and $M$ be smooth manifolds without boundary.
\begin{enumerate}[label=(\alph*)]
\item Every immersion $F\colon N\to M$ restricts to an embedding
on a sufficiently small neighborhood of each point of $N$.
\item Every injective proper immersion $F\colon N\to M$ is an
embedding. In particular, an injective immersion with compact
domain is an embedding.
\end{enumerate}
\end{proposition}
\begin{proof}
For (a), the rank theorem expresses $F$ locally as $u\mapsto(u,0)$.
Projection onto the first $\dim N$ coordinates supplies its
continuous inverse on this image.

For (b), first observe that a continuous proper map between
manifolds is closed. Let $C\subseteq N$ be closed and let
$y\in\overline{F(C)}$. Charts give a countable neighborhood basis
at $y$, so there are $x_j\in C$ with $F(x_j)\to y$. Choose a
neighborhood of $y$ with compact closure $K$. After discarding
finitely many terms, $x_j\in F^{-1}(K)$. This set is compact by
properness and is sequentially compact by compactness and first
countability. A subsequence converges to $x\in C$, and continuity
gives $F(x)=y$. Thus $F(C)$ is closed. If $F$ is also injective,
the bijection $N\to F(N)$ is continuous and closed, hence a
homeomorphism. Together with the immersion hypothesis, this proves
that it is an embedding. If $N$ is compact, the preimage of a compact
subset of $M$ is closed in $N$ and therefore compact, so $F$ is proper.
\end{proof}

\begin{proposition}[Tangent space of an embedded submanifold]
\label{prop:tangente-subvariedad-encajada-finita}
\index{tangent space!of an embedded submanifold}
Let $M$ be a smooth manifold without boundary, let $S\subseteq M$ be an
embedded submanifold without boundary of dimension $k$, and let
$\iota\colon S\hookrightarrow M$ be the inclusion. Then
$d\iota_p$ is injective and gives the identification
\[
 T_pS\cong d\iota_p(T_pS)\subseteq T_pM.
\]
In a slice chart $(x^1,\ldots,x^n)$,
\[
 d\iota_p(T_pS)
 =\operatorname{span}\left\{
 \left.\boldsymbol{\partial}_{1}\right|_p,\ldots,
 \left.\boldsymbol{\partial}_{k}\right|_p\right\}.
\]
\end{proposition}
\begin{proof}
In slice coordinates, $\iota$ is the linear inclusion
$\mathbb R^k\hookrightarrow\mathbb R^n$. Its differential has the same
form, proving both injectivity and the formula.
\end{proof}

\begin{proposition}[Tangent characterization and local extensions]
\label{prop:tangentes-extensiones-subvariedad}
Let $S^k\subseteq M$ (of dimension $n$) be an embedded submanifold, both without
boundary. Identify $T_pS$ with its image in $T_pM$.
\begin{enumerate}[label=(\alph*)]
\item A vector $v\in T_pM$ belongs to $T_pS$ if and only if it is
the velocity at zero of a smooth curve
$c\colon(-\delta,\delta)\to S$ with $c(0)=p$.
\item In a slice chart $(U,x^1,\ldots,x^n)$, the map
$h=(x^{k+1},\ldots,x^n)$ is a submersion, $U\cap S=h^{-1}(0)$,
and $T_pS=\ker dh_p$ for $p\in U\cap S$.
\item Every smooth section of $TM\restriction_S$ extends, near
each point of $S$, to a smooth vector field on $M$. In particular,
this holds for tangent vector fields on $S$.
\end{enumerate}
\end{proposition}
\begin{proof}
In slice coordinates, $T_pS$ consists of the vectors whose last
$n-k$ components vanish. Each is the velocity of the coordinate
line $t\mapsto\phi^{-1}(\phi(p)+t\,d\phi_pv)$ for sufficiently
small $|t|$. Conversely, differentiating the last components of a
curve in $S$, which vanish identically, gives the same condition.
This proves (a), and the differential of the projection $h$ proves (b).

For (c), shrink the chart image to a product box
$A\times B\subseteq\mathbb R^k\times\mathbb R^{n-k}$. Write
the section as
$\mathbf V=\displaystyle\sum_{i=1}^n V^i(u)\boldsymbol\partial_i$
on $S\cap U$. Extend each coefficient by
$\widetilde V^i(u,z):=V^i(u)$ for $(u,z)\in A\times B$.
Then $\widetilde{\mathbf V}:=
\displaystyle\sum_{i=1}^n\widetilde V^i\boldsymbol\partial_i$
is the desired extension. This is a local statement; it does not
require $S$ to be closed in $M$.
\end{proof}

\begin{proposition}[Smooth graphs]
\label{prop:grafica-subvariedad-encajada-finita}
\index{graph of a smooth map@graph of a smooth map}
Let $M$ and $N$ be smooth manifolds without boundary. If $U\subseteq M$ is
open and $f\colon U\to N$ is smooth, then
\[
 \Gamma_f:=\{(p,f(p))\mid p\in U\}\subseteq M\times N
\]
is an embedded submanifold diffeomorphic to $U$. At the point $(p,f(p))$,
\[
 T_{(p,f(p))}\Gamma_f
 =\{(v,df_pv)\mid v\in T_pM\}.
\]
\end{proposition}

\begin{proof}
The map $F\colon U\to M\times N$, $F(p)=(p,f(p))$, is injective, and
$dF_p(v)=(v,df_pv)$ is injective. The restriction of the first projection
to $\Gamma_f$ is the continuous inverse of $F$, so $F$ is an embedding. The
tangent-space formula follows from the preceding proposition.
\end{proof}

\begin{definition}\label{def:valor-regular-finito}
\index{regular value}\index{regular level set}
Let $F\colon M\to N$ be a smooth map between manifolds with or without boundary, and let
$q\in N$. We say that $q$ is a \textbf{regular value} of $F$ if
$dF_p$ is surjective for each $p\in F^{-1}(q)$. If the preimage is
empty, the condition is considered satisfied.
\end{definition}

\begin{theorem}[Regular preimage theorem]
\label{teo:preimagen-regular-finita}
\index{regular preimage theorem}
Let $M$ and $N$ be smooth manifolds without boundary, of dimensions $m$ and $n$, respectively, and let
$F\colon M\to N$ be smooth. If $q$ is a regular value, then
$F^{-1}(q)$ is either empty or a properly embedded submanifold of dimension
$m-n$, and
\[
 T_pF^{-1}(q)=\ker dF_p,\qquad p\in F^{-1}(q).
\]
\end{theorem}
\begin{proof}
At each point of the level set, Theorem~\ref{teo: teo del rango} gives
coordinates in which $F(x^1,\ldots,x^m)=(x^1,\ldots,x^n)$. The level set is
then the slice $x^1=\cdots=x^n=0$, and the preceding local criterion
gives it its embedded structure. The same coordinate expression identifies its
tangent space with the kernel of $dF_p$. Since $\{q\}$ is closed, the level set is a
closed subset of $M$; the inclusion of a closed subset is proper over
every compact set.
\end{proof}

\begin{definition}[Submanifold with boundary]
\label{def:subvariedad-encajada-con-frontera}
\index{submanifold with boundary}
A smooth manifold $S$ with boundary contained in a smooth manifold $M$
is an \textbf{embedded submanifold with boundary} if the inclusion
$S\hookrightarrow M$ is a smooth embedding. If $M$ has no boundary, a
\textbf{half-slice chart} of dimension $k$ for $S$ is a chart
of $M$ in which
\[
 S\cap U
 \longleftrightarrow
 \phi(U)\cap
 \bigl(\mathbb H^k\times\{0\}^{n-k}\bigr),
\]
where $\mathbb H^k=\mathbb R^{k-1}\times[0,\infty)$.
\end{definition}

\begin{theorem}[Local criterion with boundary]
\label{teo:criterio-media-rebanada-subvariedad}
Let $M$ be a smooth manifold without boundary and let $S\subseteq M$. Then $S$
is an embedded submanifold with boundary of dimension $k$ if and only if
each point of $S$ has a slice chart or a half-slice chart of dimension $k$,
centered at that point. The points represented on the boundary of
$\mathbb H^k$ are precisely the boundary points of the resulting manifold;
the other points are its interior points. This characterization
uniquely determines the smooth structure of $S$.
\end{theorem}
\begin{proof}
First suppose the indicated charts exist. In a slice
chart, define
\[
 \phi_S:=\operatorname{pr}_{\mathbb R^k}\circ
 \phi\restriction_{U\cap S},
\]
and in a half-slice chart use the same formula, with image in
$\mathbb H^k$. These maps are homeomorphisms onto open subsets of
$\mathbb R^k$ or relatively open subsets of $\mathbb H^k$. If two of them
overlap, the change of coordinates is obtained by restricting a change of
charts on $M$ to $\mathbb R^k\times\{0\}^{n-k}$. More precisely, if
$F$ is the ambient change of coordinates, the change on $S$ is represented by
\[
 u\longmapsto \operatorname{pr}_{\mathbb R^k}F(u,0).
\]
This expression extends smoothly to an open subset of $\mathbb R^k$, as does
its inverse. Thus the restricted charts form
a smooth atlas with boundary on $S$. A smooth change of these charts
preserves interior and boundary points, so the two alternatives at a
centered chart are exclusive. In these charts, the inclusion has the
expression $u\mapsto(u,0)$, so it is an immersion. Since $S$ carries the
subspace topology, the inclusion is also a homeomorphism onto its
image and is therefore an embedding.

Conversely, suppose the inclusion
$\iota\colon S\hookrightarrow M$ is an embedding. At an interior point,
Theorem~\ref{teo: teo del rango} applies, as in
Theorem~\ref{teo:criterio-cartas-rebanada-subvariedades}. Now let
$p\in\partial S$. Choose a boundary chart
\[
 \theta\colon W\longrightarrow O\cap\mathbb H^k
\]
of $S$, centered at $p$, where $O\subseteq\mathbb R^k$ is open, and a
chart $\psi$ of $M$ centered at $p$. After shrinking $W$ and the domain
of $\psi$, the coordinate map
\[
 f:=\psi\circ\iota\circ\theta^{-1}
\]
admits a smooth extension $\widehat f\colon\widehat O\to\mathbb R^n$ to
an open set $\widehat O\subseteq\mathbb R^k$ containing the origin. Since
$d\iota_p$ is injective, $D\widehat f(0)$ has rank $k$. Some minor of
order $k$ remains nonsingular in a neighborhood of $0$; hence
$\widehat f$ has constant rank $k$ there. The Euclidean rank
theorem gives local diffeomorphisms $\alpha$ and $\beta$, with
$\alpha(0)=0$, such that
\[
 (\beta\circ\widehat f\circ\alpha^{-1})(u)=(u,0)
\]
for $u$ in a neighborhood of $0$ in $\mathbb R^k$.

It remains to straighten the image under $\alpha$ of the half-space constituting the
original domain. In a neighborhood of $0$, define
\[
 r(u):=\bigl(\alpha^{-1}(u)\bigr)^k.
\]
Since $\alpha$ is a diffeomorphism, $dr_0\neq0$. Choose linear forms
$\ell^1,\ldots,\ell^{k-1}$ such that
$\ell^1,\ldots,\ell^{k-1},dr_0$ form a basis for
$(\mathbb R^k)^*$, and define
\[
 \eta(u):=(\ell^1(u),\ldots,\ell^{k-1}(u),r(u)).
\]
The differential of $\eta$ at $0$ is invertible. By the inverse
function theorem, there are open neighborhoods $A_\eta$ and $B_\eta$ of the
origin such that
\[
 \eta\colon A_\eta\longrightarrow B_\eta
\]
is a diffeomorphism. Shrink the domains of $\alpha$ and $\beta$ so that
the image of $\alpha$ is contained in $A_\eta$. Then
\[
 \eta\bigl(\alpha(O\cap\mathbb H^k)\bigr)
 =\eta(\alpha(O))\cap\mathbb H^k
\]
locally, since the last component of $\eta\circ\alpha$ is the last
coordinate of the original domain. After also restricting the normal
variables, choose a neighborhood $C$ of the origin in $\mathbb R^{n-k}$ and
define the diffeomorphism
\[
 \widetilde\eta\colon A_\eta\times C
 \longrightarrow B_\eta\times C,
 \qquad
 \widetilde\eta(u,z):=(\eta(u),z).
\]
Set $\alpha_1:=\eta\circ\alpha$ and
$\beta_1:=\widetilde\eta\circ\beta$, after the corresponding
restriction of the domain of $\beta$. Then
\[
\begin{aligned}
 (\beta_1\circ\widehat f\circ\alpha_1^{-1})(u)
 &=\widetilde\eta\bigl(
   \beta\circ\widehat f\circ\alpha^{-1}(\eta^{-1}(u))
   \bigr)\\
 &=\widetilde\eta(\eta^{-1}(u),0)=(u,0).
\end{aligned}
\]
Thus the ambient chart $\beta_1\circ\psi$ sends the local image of
$S$ to $\mathbb H^k\times\{0\}^{n-k}$. Since $W$ is open in $S$ and $S$
has the topology induced by $M$, there exists an ambient open set whose
intersection with $S$ is contained in $W$. Restricting the chart once more
excludes the points of $S$ outside the local image and gives
a half-slice chart centered at $p$.

The construction also proves uniqueness. In any
smooth structure with boundary for which the inclusion is an embedding,
charts obtained by restricting slice charts and half-slice
charts are smooth. Conversely, these charts generate the smooth
structure because they cover $S$. Thus every structure with the stated
property has the same maximal atlas.
\end{proof}

\begin{theorem}[The boundary as a submanifold]
\label{teo:frontera-subvariedad-encajada}
\index{boundary!as an embedded submanifold}
If $M$ is a smooth manifold with boundary of dimension $n$, then
$\partial M$ has a unique smooth manifold structure without boundary of
dimension $n-1$ for which the inclusion
$\iota_\partial\colon\partial M\hookrightarrow M$ is a proper embedding.
In a boundary chart $(x^1,\ldots,x^n)$,
\[
 d(\iota_\partial)_p(T_p\partial M)
 =\operatorname{span}\left\{
 \left.\boldsymbol{\partial}_{1}\right|_p,\ldots,
 \left.\boldsymbol{\partial}_{n-1}\right|_p\right\}.
\]
\end{theorem}
\begin{proof}
For each boundary chart $(U,\phi)$, the map
\[
 \phi_\partial
 :=\operatorname{pr}_{\mathbb R^{n-1}}\circ
   \phi\restriction_{U\cap\partial M}
\]
is a homeomorphism onto an open subset of $\mathbb R^{n-1}$. Smooth invariance
of the boundary ensures that the changes of coordinates preserve
$x^n=0$; their restrictions, and those of their inverses, are smooth. These
charts give $\partial M$ the asserted structure. In these charts, the inclusion
is $(x',0)\mapsto(x',0)$, proving the tangent-space assertion and
uniqueness. The inclusion is an embedding because we use the subspace
topology. Finally, $\partial M$ is closed in $M$, so its
intersection with any compact set is compact and the inclusion is proper.
\end{proof}

\begin{theorem}[Regular preimage in a manifold with boundary]
\label{teo:preimagen-regular-con-frontera}
Let $M$ be a smooth manifold of dimension $m$ with boundary, let $N$ be a manifold of dimension $n$ without
boundary, and let $F\colon M\to N$ be smooth. If $q\in N$ is a regular value of both
$F$ and $F\restriction_{\partial M}$, then $F^{-1}(q)$ is a
properly embedded submanifold with boundary of dimension $m-n$, and
\[
 \partial\bigl(F^{-1}(q)\bigr)
 =F^{-1}(q)\cap\partial M,\qquad
 T_pF^{-1}(q)=\ker dF_p.
\]
\end{theorem}
\begin{proof}
At interior points, Theorem~\ref{teo:preimagen-regular-finita} applies.
Let $p\in F^{-1}(q)\cap\partial M$. Take boundary coordinates
$(y^1,\dots,y^{m-1},t)$ centered at $p$, with $t\geq0$, and coordinates on
$N$ centered at $q$. Denote the coordinate representation of
$F$ by $f$. Surjectivity of
\[
 d(F\restriction_{\partial M})_p
 \colon T_p\partial M\longrightarrow T_qN
\]
implies that, after reordering the tangential variables, the matrix
$\partial(f^1,\dots,f^n)/\partial(y^1,\dots,y^n)$ is invertible at the
origin. Consider
\[
 \Psi(y,t):=
 \bigl(f(y,t),y^{n+1},\dots,y^{m-1},t\bigr).
\]
Its Jacobian at the origin is invertible. Extend $f$ smoothly to
negative values of $t$ and apply the inverse function theorem to
the extension of $\Psi$. Since the last component of $\Psi$ is $t$, the
resulting restriction is a diffeomorphism between relative neighborhoods in the
half-space and preserves its boundary. In the new coordinates
$z=\Psi(y,t)$, we have
\[
 f(z)= (z^1,\dots,z^n).
\]
Consequently, the level set $f^{-1}(0)$ is described by
$z^1=\cdots=z^n=0$ and $z^m\geq0$; it is a half-slice of dimension
$m-n$, whose boundary corresponds exactly to $t=0$.

Together with the slices obtained at interior points, these charts
give an atlas on $F^{-1}(q)$: their images are relatively open subsets
of $\mathbb H^{m-n}$, and their transition maps and inverses are
restrictions of smooth ambient coordinate changes. The subspace topology
is Hausdorff and second countable, and in these charts the inclusion
is the linear inclusion of the slice. Thus it is an embedding, by the
same atlas construction as in Theorem~\ref{teo:criterio-media-rebanada-subvariedad},
now carried out in boundary charts of $M$. In the preceding
coordinates, the tangent space to the level set is the kernel of the linear
projection representing $dF_p$, whence
$T_pF^{-1}(q)=\ker dF_p$. Finally, $F^{-1}(q)$ is closed because
$\{q\}$ is closed; thus the preimage of a compact set under the inclusion is a
closed subset of a compact set, and is therefore compact. The
inclusion is proper.
\end{proof}

Regular level sets are a special case of a condition comparing
the image of a differential with the tangent space of a submanifold.
This formulation also handles intersections.

\begin{definition}[Transversality]
\label{def:transversalidad-subvariedad}
Let $F\colon M\to P$ be a smooth map and let $S\subseteq P$ be
an embedded submanifold without boundary, with $P$ without boundary.
We say that $F$ is \textbf{transverse to $S$}, and write
$F\pitchfork S$, if
\[
dF_p(T_pM)+T_{F(p)}S=T_{F(p)}P
\qquad\text{for every }p\in F^{-1}(S).
\]
This holds vacuously if $F^{-1}(S)=\varnothing$.
\end{definition}

\begin{theorem}[Transverse preimage theorem]
\label{teo:preimagen-transversal-subvariedad}
Let $S\subseteq P$ be as above, of codimension $q$, and let
$F\colon M\to P$ be smooth, with $\dim M=m$, and transverse to $S$. If $M$ has no
boundary, $F^{-1}(S)$ is an embedded submanifold without boundary
of codimension $q$, or is empty. If $M$ has boundary, assume also
$F\restriction_{\partial M}\pitchfork S$; then the preimage is
an embedded submanifold with boundary and
\[
\partial(F^{-1}(S))=F^{-1}(S)\cap\partial M.
\]
In either case,
\[
T_pF^{-1}(S)=(dF_p)^{-1}(T_{F(p)}S),
\qquad p\in F^{-1}(S).
\]
If $S$ is closed in $P$, the inclusion of the preimage is proper.
\end{theorem}
\begin{proof}
In a slice chart around $F(p)$, write $S\cap U=h^{-1}(0)$, where
$h\colon U\to\mathbb R^q$ is projection onto the normal coordinates.
Then $\ker dh_{F(p)}=T_{F(p)}S$ and $dh_{F(p)}$ is surjective.
By linear algebra, $dh_{F(p)}\circ dF_p$ is surjective if and only
if $dF_p(T_pM)+\ker dh_{F(p)}=T_{F(p)}P$. Thus zero is a regular
value of $h\circ F$ on $F^{-1}(U)$. In the boundary case, the same
argument applied to the restriction proves regularity on $\partial M$.
Theorems~\ref{teo:preimagen-regular-finita} and
\ref{teo:preimagen-regular-con-frontera} give the local charts,
dimension, and boundary formula. Their structures agree on overlaps
because coordinates on the preimage are restrictions of ambient
coordinates; hence they define the global embedded structure.
The tangent formula follows from
$\ker d(h\circ F)_p=(dF_p)^{-1}(\ker dh_{F(p)})$.
If $S$ is closed, so is $F^{-1}(S)$; its intersection with any compact
subset of $M$ is compact, proving properness.
\end{proof}

\begin{corollary}[Transverse intersection]
\label{cor:interseccion-transversal-subvariedades}
Let $S,T$ be embedded submanifolds without boundary of a manifold
$M$ without boundary, with $T_pS+T_pT=T_pM$ at every $p\in S\cap T$.
Then $S\cap T$ is an embedded submanifold without boundary, or is
empty, and
\[
\operatorname{codim}_M(S\cap T)
=\operatorname{codim}_M S+\operatorname{codim}_M T,
\qquad T_p(S\cap T)=T_pS\cap T_pT.
\]
\end{corollary}
\begin{proof}
Apply the theorem to the inclusion $S\hookrightarrow M$ and the
submanifold $T$. The resulting dimension is
$\dim S-\operatorname{codim}_M T$, and the tangent formula gives
the asserted intersection. The composition of embeddings
$S\cap T\hookrightarrow S\hookrightarrow M$ is an embedding:
the induced topologies agree and the differentials are injective.
\end{proof}

The boundary hypothesis cannot be omitted. On
$M=\mathbb R\times[0,\infty)$, the map $F(x,t)=t$ has zero as a
regular value, but $F\restriction_{\partial M}$ is constant.
The level $F^{-1}(0)=\mathbb R\times\{0\}$ is a line without
boundary: its intrinsic boundary is empty, whereas its intersection
with $\partial M$ is the whole line. More generally, the boundary of
a manifold and the topological boundary of a subset are distinct
notions. For example, a circle embedded in $\mathbb R^2$ has no
intrinsic boundary, although it is its own topological boundary.

The double replaces a manifold with boundary by a manifold without
boundary containing two copies of the original manifold. This construction is
useful, for example, for extending smooth objects across the boundary and
locally applying results formulated for manifolds without boundary.

\begin{theorem}[Collar and smooth double]
\label{teo:doble-suave-extension-frontera}
\index{collar of a manifold}\index{double of a manifold with boundary}
Let $M$ be a smooth manifold with boundary.
\begin{enumerate}[label=(\alph*)]
\item There exist $\varepsilon>0$ and a smooth embedding
\[
 c\colon\partial M\times[0,\varepsilon)\longrightarrow M,\qquad
 c(x,0)=x,
\]
whose image is a neighborhood of $\partial M$ in $M$.
\item If $M_+$ and $M_-$ are two copies of $M$, the quotient
\[
 DM:=(M_+\sqcup M_-)\big/\bigl(x_+\sim x_-\text{ for }
 x\in\partial M\bigr)
\]
admits a smooth manifold structure without boundary. The inclusions of
both copies are smooth embeddings with boundary, and their boundaries have the
same image, an embedded hypersurface of $DM$.
\end{enumerate}
The diffeomorphism type of $DM$ does not depend on the chosen collar.
\end{theorem}

\begin{proof}
Part (a) is the collar theorem
\cite[Theorem~9.25]{LeeS}. For part (b), apply the smooth manifold gluing
theorem to the two copies of $M$ and the identity diffeomorphism
of $\partial M$; see
\cite[Theorem~9.29 and Example~9.32]{LeeS}. These results construct the
smooth structure on the quotient, prove the assertions about the
inclusions, and show that two choices of collar produce diffeomorphic
doubles.
\end{proof}

\begin{corollary}[Smooth extension to the double]
\label{cor:extension-suave-al-doble}
Let $M$ be a smooth manifold with boundary and let $DM$ be its double. Every
smooth function and every smooth tensor field on one copy of $M$ admit
a smooth extension to $DM$. If $\mathbf g$ is a Riemannian metric on
$M$, a symmetric tensor extension $\widetilde{\mathbf
g}$ can be chosen that is Riemannian on an open neighborhood of that copy.
\end{corollary}

\begin{proof}
Identify $M$ with the copy $M_+\subset DM$ and let $\mathbf T$ be a smooth tensor
field on $M_+$. By the definition of smoothness up to the
boundary, for every $p\in M_+$ there exist an open set $V_p\subset DM$
containing $p$ and a smooth tensor field $\widetilde{\mathbf T}_p$ on
$V_p$ whose restriction to $V_p\cap M_+$ agrees with $\mathbf T$. At
interior points, we simply take $V_p\subset \operatorname{Int}(M_+)$; at
boundary points, we extend the local components of
$\mathbf T$ to negative values of the normal coordinate in a chart of the
double.

The family
\[
 \{V_p:p\in M_+\}\cup\{DM\setminus M_+\}
\]
is an open cover of $DM$. Let
$(\chi_p)_{p\in M_+}\cup\{\chi_0\}$ be a smooth, locally finite partition of unity
subordinate to this cover. Each product
$\chi_p\widetilde{\mathbf T}_p$, extended by zero outside $V_p$, is
a smooth tensor field on $DM$. By local finiteness, the sum
\[
 \widetilde{\mathbf T}:=
 \displaystyle\sum_{p\in M_+}\chi_p\widetilde{\mathbf T}_p
\]
is well defined and smooth. If $x\in M_+$, then $\chi_0(x)=0$, and
each $\widetilde{\mathbf T}_p(x)$ occurring in the sum agrees with
$\mathbf T(x)$; consequently,
\[
 \widetilde{\mathbf T}(x)
 =\left(\displaystyle\sum_{p\in M_+}\chi_p(x)\right)\mathbf T(x)
 =\mathbf T(x).
\]
The same argument, with tensors of type $(0,0)$, treats functions.

If $\mathbf T=\mathbf g$, replace the resulting extension by its symmetric
part. The extension still agrees with $\mathbf g$ on $M_+$.
Since the set of positive definite symmetric bilinear forms is
open in each fiber and the extension is continuous, the set of points
of $DM$ where $\widetilde{\mathbf g}$ is positive definite is an open set
containing $M_+$. Therefore, $\widetilde{\mathbf g}$ is a
Riemannian metric on an open neighborhood of that copy.
\end{proof}

 \section{Vector bundles and the cotangent bundle}

The tangent bundle suggests a more general construction. Over each point of a manifold we may place a vector space, provided these fibers can be identified locally in a smooth manner and the identifications are linear on each fiber. This leads to vector bundles, the natural setting for fields, differential forms, tensors, spinors, and the geometric unknowns that will appear later.

 \begin{definition}\label{def:nociones-fundamentales-haz-vectorial-real-de-rango-sobre}\index{real vector bundle of rank}
 Let $M$ be a topological space. A \textbf{real vector bundle of rank $k$ over $M$} is a topological space $\mathbf{E}$ together with a continuous surjective function $\pi\colon \mathbf{E}\longrightarrow M$ such that:
 \begin{enumerate}[label=(\alph*)]
 \item For each $p\in M$, the fiber $\mathbf{E}_{p}=\pi^{-1}(p)$ over $p$ is equipped with the structure of a real vector space of dimension $k$.
 \item For each $p\in M$, there exist a neighborhood $U$ of $p$ in $M$ and a homeomorphism $\Phi\colon \pi^{-1}(U)\longrightarrow U\times \mathbb{R}^{k}$ such that $\pi_{U}\circ \Phi=\pi$, where $\pi_{U}\colon U\times \mathbb{R}^{k}\longrightarrow U$ is the projection, and for each $q\in U$, the restriction of $\Phi$ to $\mathbf{E}_{q}$ is an isomorphism from $\mathbf{E}_{q}$ to $\{q\}\times\mathbb{R}^{k}\cong \mathbb{R}^{k}$. ($\Phi$ is called a \textbf{local trivialization of $\mathbf{E}$ over $U$.})
 \end{enumerate}
 \end{definition}
 \begin{definition}\label{def:nociones-fundamentales-haz-vectorial-suave-de-rango-sobre}\index{smooth vector bundle of rank}\index{vector bundle!smooth}\glsadd{haz-vectorial}
 Let $M$ and $\mathbf{E}$ be smooth manifolds with or without boundary such that $\pi\colon \mathbf{E}\longrightarrow M$ is a continuous surjective map making $\mathbf{E}$ a real vector bundle of rank $k$ over $M$. We say that $\mathbf{E}$ is a \textbf{smooth vector bundle of rank $k$ over $M$} if
 \begin{enumerate}[label=(\alph*)]
 \item $\pi$ is a smooth map.
 \item For each $p\in M$, there exist a neighborhood $U$ of $p$ and a local trivialization $\Phi\colon \pi^{-1}(U)\longrightarrow U\times\mathbb{R}^{k}$ of $\mathbf{E}$ over $U$ such that $\Phi$ is a diffeomorphism. This map is called a \textbf{smooth local trivialization of $\mathbf{E}$ over $U$}.
 \end{enumerate}
 \end{definition}

\begin{definition}[Smooth complex vector bundle]
\label{def:haz-vectorial-complejo-suave}
Let $M$ be a smooth manifold with or without boundary. A \emph{smooth complex
vector bundle of rank $r$} over $M$ is a smooth map
$\pi\colon\mathbf E\to M$ admitting, around each point of $M$, a
smooth trivialization
\[
 \Phi\colon\pi^{-1}(U)\longrightarrow U\times\mathbb C^r
\]
whose restriction to each fiber is complex linear. In particular, the
total space is regarded as a smooth real manifold, and the transition
functions take values in $GL(r,\mathbb C)$.
\end{definition}

\begin{figura}[h]
    \includegraphics[scale=1]{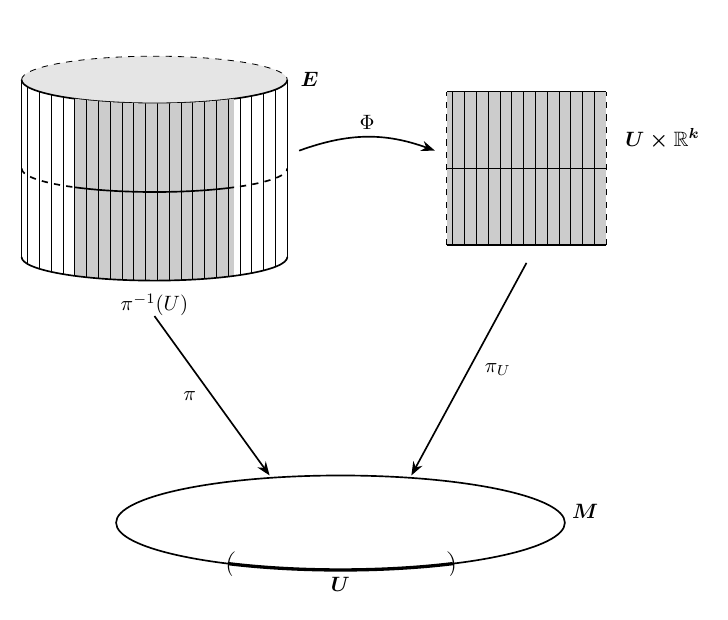}
    \caption{Local trivialization of a vector bundle.}
    \label{fig:trivializacion-local-haz-vectorial}
\end{figura}
\begin{remark}\label{obs:nociones-fundamentales-ocasiones-no-sobrecargar-notacion-siempre-no}
When no ambiguity arises, we omit the projection when referring to a vector bundle, and may even omit the base space $M$. Thus, if $\pi\colon \mathbf{E}\longrightarrow M$ is a vector bundle, we write $\mathbf{E}\longrightarrow M$ or $(\mathbf{E},M)$ or simply $\mathbf{E}$.
\end{remark}

The simplest example of a vector bundle is what is known as the \textit{trivial bundle of rank $k$ over $M$}.
\begin{example}\label{ej:nociones-fundamentales-haz-trivial-de-rango-sobre}
 If $M$ is a smooth manifold with or without boundary, the bundle $\underline{\mathbb{R}}^{k}_{M}=M\times \mathbb{R}^{k}$ is a smooth vector bundle of rank $k$ over $M$. It is called the \textit{\textbf{trivial bundle of rank $k$ over $M$}}.
\end{example}
The vector bundle in this example is trivial in the sense that it admits a global trivialization (a local trivialization with $U=M$).

The tangent bundle provides another example of a vector bundle.
 \begin{proposition}\label{prop:nociones-fundamentales-haz-tangente-haz-vectorial-suave-rango}
 Let $M$ be a smooth manifold with or without boundary of dimension $n$. Its tangent
 bundle $TM$ is a smooth vector bundle of rank $n$ over $M$.
 \end{proposition}

 \begin{proof}
 Given a coordinate chart $(U,\phi)$ of $M$ such that $\phi=(x^{1},\dots,x^{n})$ and $p\in U$, define the map $\Phi\colon \pi^{-1}(U)\longrightarrow U\times \mathbb{R}^{n}$ by \[\Phi\left( p,v^{i}\boldsymbol{\partial}_{i}\biggr|_{p}\right)=(p,(v^{1},\dots,v^{n})).\]
 The map $\Phi$ is linear on the fibers of $\pi$ and satisfies $\pi_{U}\circ \Phi=\pi$.
 If $(U,\phi)$ is an interior chart, the tangent-bundle chart in
 Proposition~\ref{estructurasuavedelhaztangente} satisfies
 \[
 \widetilde{\phi}=(\phi\times\text{Id}_{\mathbb R^n})\circ\Phi.
 \]
 If it is a boundary chart, the tangent-bundle chart places the
 vector coordinates first, followed by those of the base point. Therefore,
 \[
 \widetilde{\phi}=S\circ(\phi\times\text{Id}_{\mathbb R^n})\circ\Phi,
 \qquad S(x,v)=(v,x).
 \]
 The interchange $S$ is a linear isomorphism of $\mathbb R^{2n}$ and sends
 $\phi(U)\times\mathbb R^n$ diffeomorphically onto
 $\mathbb R^n\times\phi(U)$, a relatively open subset of $\mathbb H^{2n}$.
 In both cases, $\Phi$ is a composition of diffeomorphisms, and hence
 is a smooth local trivialization.
 \end{proof}
 Every nontrivial vector bundle requires more than one local trivialization. The following lemma shows that the composition of two smooth local trivializations has a simple form where they overlap.

 \begin{lemma}\label{lema existnecia mapeo GL}
Let $M$ be a smooth manifold with or without boundary and let
$\pi\colon \mathbf{E}\longrightarrow M$ be a smooth vector bundle of rank
$k$. Suppose $\Phi\colon \pi^{-1}(U)\longrightarrow U\times
\mathbb{R}^{k}$ and $\Psi\colon \pi^{-1}(V)\longrightarrow
V\times\mathbb{R}^{k}$ are smooth local trivializations of $\mathbf{E}$
with $U\cap V\neq \varnothing$. Then there exists a smooth map
$\tau\colon U\cap V\longrightarrow GL(k,\mathbb{R})$ such that
\[
\Phi\circ \Psi^{-1}\colon (U\cap V)\times \mathbb{R}^{k}\longrightarrow (U\cap V)\times \mathbb{R}^{k}
\]
has the form
\[
(\Phi\circ \Psi^{-1})(p,v)=(p,\tau(p)v),
\]
where $\tau(p)v$ denotes the usual action of the matrix $\tau(p)$ in $k\times k$ on $v\in \mathbb{R}^{k}$.
\end{lemma}

\begin{proof}
The following diagram commutes:
\[
\begin{tikzcd}
(U \cap V)\times \mathbb{R}^k \arrow[dr, "\pi_1"'] &
\pi^{-1}(U \cap V) \arrow[l, "\Psi"'] \arrow[r, "\Phi"] \arrow[d, "\pi"] &
(U \cap V)\times \mathbb{R}^k \arrow[dl, "\pi_1"] \\
& U \cap V &
\end{tikzcd}
\]
where the maps at the top are interpreted as the restrictions of $\Psi$ and $\Phi$ to $\pi^{-1}(U\cap V)$. Commutativity of the diagram tells us that
\[
\pi_{1}\circ (\Phi\circ \Psi^{-1})=\pi_{1},
\]
where $\pi_{1}$ is projection onto the first coordinate.

This implies that there exists a map $\sigma\colon (U\cap V)\times \mathbb{R}^{k}\longrightarrow \mathbb{R}^{k}$ such that
\[
(\Phi\circ \Psi^{-1})(p,v)=(p,\sigma(p,v)).
\]
By construction, $\sigma$ is the second coordinate function of $\Phi\circ \Psi^{-1}$; that is,
\[
\sigma=\pi_{2}\circ (\Phi\circ \Psi^{-1}),
\]
where $\pi_{2}$ is projection onto the second coordinate.

For each $p\in U\cap V$, the map $v\mapsto \sigma(p,v)$ from $\mathbb{R}^{k}$ to itself is linear and invertible, since $\Phi$ and $\Psi$ are local trivializations; hence $\Phi\circ\Psi^{-1}$ is a fiberwise diffeomorphism preserving the first coordinate. In particular, there is a unique matrix $\tau(p)\in GL(k,\mathbb{R})$ such that
\[
\sigma(p,v)=\tau(p)v\qquad\text{for every }v\in\mathbb{R}^{k}.
\]
We now show that $p\mapsto \tau(p)$ is smooth.

Fix the standard basis $(e_{1},\dots,e_{k})$ of $\mathbb{R}^{k}$. Since $v\mapsto \sigma(p,v)$ is linear, it is determined by its values on the basis vectors, and we define $\tau(p)$ as the matrix with these values as columns:
\[
\tau(p)e_{j}=\sigma(p,e_{j})\qquad\text{for }j\in\{1,\dots,k\}.
\]
The map
\[
p\longmapsto \sigma(p,e_{j})=\pi_{2}\big((\Phi\circ \Psi^{-1})(p,e_{j})\big)
\]
is a composition of smooth functions (since $(p,e_{j})$ is constant in the second coordinate and $\Phi\circ \Psi^{-1}$ is smooth), and is therefore a smooth map from $U\cap V$ to $ \mathbb{R}^{k}$ for each $j\in\{1,\dots,k\}$.

Denoting projection onto the $i$th coordinate by $\mathrm{pr}_{i}\colon \mathbb{R}^{k}\longrightarrow\mathbb{R}$, the entries of the matrix $\tau(p)=(\tau_{ij}(p))$ are given by
\[
\tau_{ij}(p)=\mathrm{pr}_{i}\big(\sigma(p,e_{j})\big),\qquad i,j\in\{1,\dots,k\}.
\]
Each function $\tau_{ij}$ is smooth on $p$, so $p\mapsto \tau(p)$ is a smooth map.

Consequently,
\[
(\Phi\circ \Psi^{-1})(p,v)=(p,\tau(p)v)
\]
with $\tau\colon U\cap V\longrightarrow GL(k,\mathbb{R})$ a smooth function, as desired.
\end{proof}

 The following result allows us to construct smooth vector bundles.
 \begin{lemma}[Smooth vector bundle lemma]\label{lema del haz vectorial suave}\index{smooth vector bundle lemma}
 Let $\mathbb K\in\{\mathbb R,\mathbb C\}$, let $M$ be a smooth manifold
 with or without boundary, and for each $p\in M$ let $\mathbf{E}_p$ be a vector space
 of fixed dimension $r$ over $\mathbb K$. Set
 $\mathbf{E}=\coprod_{p\in M}\mathbf{E}_p$ and let $\pi\colon \mathbf{E}\longrightarrow M$ be the map
 sending each element of $\mathbf{E}_p$ to $p$. Suppose we are given:
 \begin{enumerate}[label=(\alph*)]
 \item An open cover $\{U_{\alpha}\}_{\alpha\in A}$ of $M$.
 \item For each $\alpha\in A$, a fiber-preserving bijection
 \[
 \Phi_\alpha\colon\pi^{-1}(U_\alpha)
 \longrightarrow U_\alpha\times\mathbb K^r
 \]
 whose restriction to $\mathbf{E}_p$ is a $\mathbb K$-linear isomorphism from
 $\mathbf{E}_p$ onto $\{p\}\times\mathbb K^r$.
 \item For $U_\alpha\cap U_\beta\neq\varnothing$, a smooth map
 $g_{\alpha\beta}\colon U_\alpha\cap U_\beta
 \longrightarrow GL(r,\mathbb K)$ such that
 \[
 (\Phi_\alpha\circ\Phi_\beta^{-1})(p,z)
 =(p,g_{\alpha\beta}(p)z).
 \]
 \end{enumerate}

 Then $\mathbf{E}$ has a unique topology and a unique real smooth structure
 making it a smooth manifold with or without boundary, for
 which $\pi$ is smooth and the maps $\Phi_\alpha$ are smooth,
 $\mathbb K$-linear local trivializations. Addition
 $\mathbf{E}\times_M \mathbf{E}\to \mathbf{E}$ and scalar multiplication
 $\mathbb K\times \mathbf{E}\to \mathbf{E}$ are smooth. In particular, $\mathbf{E}\to M$ is a smooth
 real vector bundle if $\mathbb K=\mathbb R$ and a smooth complex vector bundle if
 $\mathbb K=\mathbb C$.
\end{lemma}
\begin{proof}
Let $d=\dim_{\mathbb R}\mathbb K$. For each $\alpha$ and each chart
$(V,\varphi)$ of $M$ with $V\subseteq U_\alpha$, define
\[
 \widetilde\varphi_{\alpha,V}
 :=(\varphi\times\operatorname{Id}_{\mathbb K^r})
 \circ\Phi_\alpha
 \colon\pi^{-1}(V)\longrightarrow
 \varphi(V)\times\mathbb K^r.
\]
Identify $\mathbb K^r$ with $\mathbb R^{dr}$. If
$(V,\varphi)$ is an interior chart, its image is open in
$\mathbb R^{n+dr}$; if it is a boundary chart, reordering the
coordinates identifies
$\varphi(V)\times\mathbb R^{dr}$ with a relatively open subset of
$\mathbb H^{n+dr}$.

Consider two charts $(V,\varphi)$ and $(W,\psi)$ subordinate,
respectively, to $U_\alpha$ and $U_\beta$. On the intersection, the change
of coordinates from the second chart to the first is
\[
 (y,z)\longmapsto
 \left(
   \varphi\circ\psi^{-1}(y),
   g_{\alpha\beta}(\psi^{-1}(y))z
 \right).
\]
It is smooth as a real map, and its inverse is obtained by interchanging
$\alpha,\beta$; indeed, the given bijections imply
$g_{\beta\alpha}=g_{\alpha\beta}^{-1}$. The images of the intersections
are open because they have the form
$\varphi(V\cap W)\times\mathbb K^r$.

This family of charts on the total space has a countable subfamily
covering $\mathbf{E}$: simply choose a countable subcover of the domains $V$ on the
base, since $M$ is second countable. Two points of $\mathbf{E}$ lying in the same
fiber belong to a common chart $\pi^{-1}(V)$. If their base points are
distinct, first choose disjoint coordinate domains in $M$,
restricted to their respective sets $U_\alpha$, and then take their preimages.
This verifies the second countability and separation assumptions of
Lemma~\ref{lema de la carta suave}. The lemma provides a Hausdorff,
second countable topology and a smooth structure on $\mathbf{E}$ for which
the maps $\widetilde\varphi_{\alpha,V}$ are charts.

Each $\Phi_\alpha$ is a diffeomorphism: after restriction to
$\pi^{-1}(V)$, with $V\subseteq U_\alpha$, its expression in the newly
constructed charts and the product chart of $V\times\mathbb K^r$ is the identity.
In the same
coordinates, $\pi$ is $(x,z)\mapsto x$, and hence is smooth. Fiberwise addition and
scalar multiplication have the expressions
\[
 (x,z,w)\longmapsto(x,z+w),
 \qquad
 (\lambda,x,z)\longmapsto(x,\lambda z),
\]
which are smooth as real maps and $\mathbb K$-linear in the fiber
variables.

If another topology and smooth structure had the properties in the
statement, all the maps $\Phi_\alpha$ would be diffeomorphisms. Consequently,
the sets $\Phi_\alpha^{-1}(O)$, with $O$ open in
$U_\alpha\times\mathbb K^r$, would necessarily form a basis for the same
topology, and the charts $\widetilde\varphi_{\alpha,V}$ would belong to the same
smooth structure. This proves uniqueness. In the case with boundary, the
construction also shows that
$\partial \mathbf{E}=\pi^{-1}(\partial M)$.

In the complex case, the same construction is carried out on the underlying
real manifold and preserves the complex linearity of the transition maps.
\end{proof}
 This lemma produces nontrivial vector bundles, as the following example shows:

 \begin{example}[Möbius bundle]\label{haz:mobius}
Define on $\mathbb{R}^{2}$ the relation
\[
(x,y)\sim(x',y') \Longleftrightarrow (x',y')=(x+n,(-1)^{n}y),\quad n\in\mathbb{Z}.
\]
Let $\mathbf{E}=\mathbb{R}^{2}/\sim$ and let
$q\colon\mathbb{R}^{2}\longrightarrow \mathbf{E}$ be the quotient projection.
Geometrically, $\mathbf{E}$ is obtained from $[0,1]\times\mathbb{R}$ by identifying
$(0,y)$ with $(1,-y)$; the subset
$q([0,1]\times[-r,r])$ is the compact Möbius strip of width $2r$.

Let $\pi_{1}(x,y)=x$ and $\varepsilon(x)=e^{2\pi i x}$. If $U$ is a proper open arc
of $\mathbb S^1$ admitting a continuous branch of the
argument, each component $\widetilde U$ of $\varepsilon^{-1}(U)$ is an
open interval and
$\varepsilon\restriction_{\widetilde U}\colon\widetilde U\to U$ is a
diffeomorphism; its inverse is obtained by dividing that branch of the
argument by $2\pi$. Since
$\varepsilon(x+n)=\varepsilon(x)$, the map
$\varepsilon\circ\pi_1$ is constant on the equivalence classes of $\sim$. Therefore,
\[
\pi\colon \mathbf{E}\longrightarrow\mathbb S^1,
\qquad
\pi([(x,y)])=\varepsilon(x),
\]
is well defined and satisfies $\pi\circ q=\varepsilon\circ\pi_1$. This is
summarized by the diagram:

\begin{tikzcd}
\mathbb{R}^{2} \arrow[r, "q"] \arrow[d, "\pi_{1}"']
 & \mathbf{E}=\mathbb{R}^{2}/\sim \arrow[d, "\pi", dashed] \\
\mathbb{R} \arrow[r, "\varepsilon"]
 & \mathbb{S}^{1}
\end{tikzcd}

Let \(p\in\mathbb{S}^{1}\) and choose \(x_{0}\in\mathbb{R}\) such that \(\varepsilon(x_{0})=p\). First characterize the real numbers \(x\) with \(\varepsilon(x)=p\). Since \(\varepsilon(x)=e^{2\pi i x}\), the equality \(\varepsilon(x)=\varepsilon(x_{0})\) is equivalent to \(e^{2\pi i(x-x_{0})}=1\). The real solutions of \(e^{2\pi i t}=1\) are precisely the integers \(t\in\mathbb{Z}\), so \(\varepsilon(x)=\varepsilon(x_{0})\) if and only if \(x-x_{0}\in\mathbb{Z}\), that is, if and only if \(x=x_{0}+n\) for some \(n\in\mathbb{Z}\). Consequently,
\[
\pi^{-1}(p)=\{[(x,y)]\in \mathbf{E}\mid \varepsilon(x)=p\}=\{[(x_{0}+n,y)]\mid n\in\mathbb{Z},\ y\in\mathbb{R}\}.
\]
We can rewrite this fiber by recalling the definition of $\sim$. We know that $(x_{0}+n,y)\sim (x_{0},(-1)^{n}y)$ for any $n\in \mathbb{Z}$ and $y\in \mathbb{R}$. This allows us to conclude that \[\pi^{-1}(p)=\{[(x_{0},(-1)^{n}y)]\mid n\in \mathbb{Z}, y\in \mathbb{R}\},\]

which can be rewritten as $\pi^{-1}(p)=\{[(x_{0},t)]\mid t\in \mathbb{R}\}$.

Define
$j_{p}\colon \pi^{-1}(p)\longrightarrow \mathbb{R}$
by $j_{p}([(x_{0},t)])=t$.

 We show that this does not depend on the chosen $t$: if $(x_{0},t)\sim (x_{0},t')$, there exists $n\in \mathbb{Z}$ such that $(x_{0},t')=(x_{0}+n,(-1)^{n}t)$, which gives $n=0$ and therefore $t=t'$.

The function $j_p$ is injective because
$j_p([(x_0,t)])=j_p([(x_0,s)])$ implies $t=s$, and surjective because
$a=j_p([(x_0,a)])$ for every $a\in\mathbb R$. Equip $\pi^{-1}(p)$ with
the vector space structure given by
\[[(x_{0},t)]+[(x_{0},s)]:=[(x_{0},t+s)]\]
and
\[\lambda [(x_{0},t)]=[(x_{0},\lambda t)].\]
These operations are well defined by the uniqueness of the representative of
the form $(x_0,t)$ proved above. With these operations, $j_p$ preserves addition and
scalar multiplication, and is therefore a linear isomorphism.
If another lift $x'_0=x_0+m$ is chosen, with $m\in\mathbb Z$, then
\[
[(x_0,t)]=[(x'_0,(-1)^m t)].
\]
The corresponding isomorphism $j'_p$ satisfies $j'_p=(-1)^m j_p$.
Since multiplication by $(-1)^m$ is a linear isomorphism of
$\mathbb R$, the operations of addition and scalar multiplication
transported to the fiber are unchanged. Its vector space structure
is therefore independent of the chosen lift.

Take the cover
\[
U_\pm=\mathbb{S}^{1}\setminus\{(\pm 1,0)\}.
\]
Fix components
\[
\widetilde{U}_{-}=(-\displaystyle\frac{1}{2},\displaystyle\frac{1}{2})\subset\varepsilon^{-1}(U_{-}),\qquad
\widetilde{U}_{+}=(0,1)\subset\varepsilon^{-1}(U_{+}),
\]
so that $\varepsilon|_{\widetilde{U}_{\pm}}\colon \widetilde{U}_{\pm}\longrightarrow U_{\pm}$ are homeomorphisms. Define the trivializations
\[
\Phi_{\pm}\colon \pi^{-1}(U_{\pm})\longrightarrow U_{\pm}\times\mathbb{R},\quad [(x,y)]\mapsto(\varepsilon(x),y),\ \ x\in\widetilde{U}_{\pm},
\]

Both satisfy $\pi=\pi_{1}\circ\Phi_{\pm}$. For $p\in U_\pm$, let
$x_\pm\in\widetilde U_\pm$ be the unique point with $\varepsilon(x_\pm)=p$.
If $x_0\in\mathbb R$ is any other lift of $p$, there exists
$m\in\mathbb Z$ such that $x_0=x_\pm+m$, and then
\[
\Phi_\pm([(x_0,t)])
=\Phi_\pm([(x_\pm,(-1)^m t)])
=(p,(-1)^m t).
\]
In particular, using the representative with $x_0=x_\pm$, the second
component is $t$. With any other lift, it changes only
by the linear isomorphism $t\mapsto(-1)^m t$. Thus the restriction
of $\Phi_\pm$ to each fiber is a vector space isomorphism.
\begin{center}
\begin{tikzcd}[row sep=1.5em, column sep=3em]
& \pi^{-1}(U_{\pm}) \arrow[dl, "\pi"'] \arrow[dr, "\Phi_{\pm}"] & \\
U_{\pm} & & U_{\pm}\times\mathbb{R} \arrow[ll, "\pi_{1}"]
\end{tikzcd}
\end{center}

We next compute the transition maps. Let $(p,v)\in (U_{+}\cap U_{-})\times\mathbb{R}$. To invert $\Phi_{-}$, choose the representative
\[
\Phi_{-}^{-1}(p,v)=[(x_{-},v)]\quad\text{where}\quad x_{-}=(\varepsilon|_{\widetilde{U}_{-}})^{-1}(p)\in (-\displaystyle\frac{1}{2},\displaystyle\frac{1}{2})\setminus\{0\}.
\]
To then apply $\Phi_{+}$, express the same class using a representative whose first component lies in $\widetilde{U}_{+}=(0,1)$:
\[
x_{+}=
\begin{cases}
x_{-}+1, & \text{if }x_{-}\in(-\displaystyle\frac{1}{2},0),\\
x_{-}, & \text{if }x_{-}\in(0,\displaystyle\frac{1}{2}),
\end{cases}
\qquad
[(x_{-},v)]=
\begin{cases}
[(x_{+},-v)], &\text{if } x_{-}\in (-\frac{1}{2},0),\\
[(x_{+},v)], &\text{if } x_{-}\in(0,\frac{1}{2}),
\end{cases}
\]
since $(x_{-},v)\sim(x_{-}+1,-v)$ by the equivalence relation. Consequently,
\[
(\Phi_{+}\circ \Phi_{-}^{-1})(p,v)=
\begin{cases}
(p,-v), & \text{if }x_{-}\in(-\displaystyle\frac{1}{2},0),\\
(p,v), & \text{if }x_{-}\in(0,\displaystyle\frac{1}{2}).
\end{cases}
\]
This can be written as
\[
(\Phi_{+}\circ \Phi_{-}^{-1})(p,v)=(p,\tau_{+-}(p)v),
\qquad
\tau_{+-}\colon U_{+}\cap U_{-}\longrightarrow GL(1,\mathbb{R}),
\]
where $\tau_{+-}(p)=-1$ on the component corresponding to $x_{-}\in(-\displaystyle\frac{1}{2},0)$ and $\tau_{+-}(p)=1$ on the component corresponding to $x_{-}\in(0,\displaystyle\frac{1}{2})$. The image of $\tau_{+-}$ is the subset $\{+1,-1\}\subset GL(1,\mathbb{R})$, and $\tau_{+-}$ is smooth because it is constant on each connected component of $U_{+}\cap U_{-}$.

These trivializations and transition maps satisfy the assumptions of
Lemma~\ref{lema del haz vectorial suave}. Therefore,
$(\mathbf{E},\pi,\mathbb{S}^{1})$ is a smooth vector bundle of rank $1$.
To verify that it is nontrivial, suppose it admits a nowhere vanishing
continuous section. Lifting it through $q$ represents the section
by a continuous function $a\colon\mathbb R\to\mathbb R$
satisfying
\[
 a(x+1)=-a(x),
\]
because of the relation $(x,y)\sim(x+1,-y)$. If $a(0)>0$, then
$a(1)<0$, and if $a(0)<0$, then $a(1)>0$; in both cases, the intermediate
value theorem gives $x\in[0,1]$ with $a(x)=0$. The case $a(0)=0$
is immediate. Thus every continuous section vanishes somewhere. A
trivial rank-one bundle has the nowhere vanishing global section $p\mapsto(p,1)$;
therefore, $\mathbf E$ is nontrivial.
\end{example}
\begin{figura}[H]
    \includegraphics[scale=0.8]{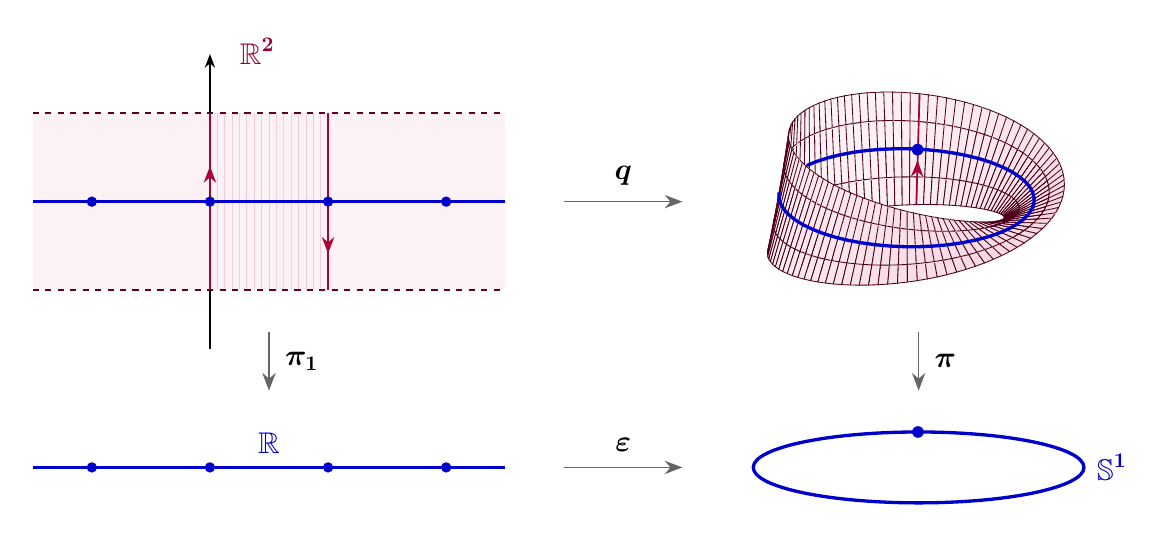}
    \caption{The Möbius bundle and $[0,1]\times\mathbb{R}$ as a fundamental domain representing it.}
    \label{fig:haz-mobius-dominio-fundamental}
\end{figura}
Operations on vector spaces also produce new bundles, as follows. \begin{definition}\label{suma de whitney}\index{Whitney sum@Whitney sum}
 Given a smooth manifold $M$ with or without boundary and smooth vector bundles $\pi'\colon \mathbf{E}'\longrightarrow M$ and $\pi''\colon \mathbf{E}''\longrightarrow M$ of ranks $k'$ and $k''$, respectively, define the \textbf{Whitney sum of $\mathbf{E}'$ and $\mathbf{E}''$} by
 $\mathbf{E}'\oplus_W \mathbf{E}'':=\displaystyle\coprod _{p\in M}(\mathbf{E}_{p}'\oplus \mathbf{E}_{p}'')$, with projection $\pi\colon \mathbf{E}'\oplus \mathbf{E}''\longrightarrow M$ given by $\pi(p,(v',v''))=p$ for each $p\in M$.

 \end{definition}
 \begin{proposition}\label{prop:nociones-fundamentales-suma-whitney-dos-haces-vectoriales-suaves}
 Let $M$ be a smooth manifold with or without boundary. The Whitney sum of two smooth vector bundles over $M$ of ranks $k'$ and $k''$ is a smooth vector bundle over $M$ of rank $k'+k''$.
 \end{proposition}
 \begin{proof}
 For each $p\in M$, consider an open neighborhood $U$ of $p$ over which there exist local trivializations $(U,\Phi')$ of $\mathbf{E}'$ and $(U,\Phi'')$ of $\mathbf{E}''$, and define $\Phi\colon \pi^{-1}(U)\longrightarrow U\times \mathbb{R}^{k'+k''}$ by $\Phi(v',v'')=(\pi'(v'),(\pi_{R^{k'}}\circ \Phi'(v'),\pi_{\mathbb{R}^{k''}}\circ \Phi''(v''))).$

 Suppose we are given another pair of local trivializations $(\widetilde{U},\widetilde{\Phi}')$ and $(\widetilde{U},\widetilde{\Phi}'')$. Let $\tau'\colon (U\cap \widetilde{U})\longrightarrow GL(k',\mathbb{R})$ and $\tau''\colon (U\cap \widetilde{U})\longrightarrow GL(k'',\mathbb{R})$ be their respective transition functions. Then the transition function for $\mathbf{E}'\oplus_W \mathbf{E}''$ is given by:
 \[(\widetilde{\Phi}\circ \Phi^{-1})(p,(v',v''))=(p,\tau(p)(v',v'')),\] where $\tau(p)=\begin{bmatrix}
 \tau'(p) & 0 \\
 0 & \tau''(p)
 \end{bmatrix},$ is smooth. Thus Lemma~\ref{lema del haz vectorial suave} ensures that $\mathbf{E}'\oplus \mathbf{E}''$ is a vector bundle of rank $k'+k''$ over $M$.
 \end{proof}
 We now introduce \textit{local and global sections} and \textit{local and global frames}. These will allow us to describe fully the vector bundles commonly used on a smooth manifold $M$. They will also make it easier to define vector fields and tensor fields.
 \begin{definition}\label{def:nociones-fundamentales-seccion-local-de}\index{local section@local section}
 Let $\pi\colon \mathbf{E}\longrightarrow M$ be a vector bundle. A \textbf{local section of $\mathbf{E}$} is a continuous map $\boldsymbol{\sigma}\colon U\longrightarrow \mathbf{E}$ defined on an open set $U\subseteq M$ such that $\pi\circ \boldsymbol{\sigma}=\text{Id}_{U}$. If $U=M$, we also call $\boldsymbol{\sigma}$ a \textbf{global section of $\mathbf{E}$}.

 \end{definition}
 \begin{definition}\label{def:nociones-fundamentales-seccion-local-o-global-suave-de}\index{smooth local or global section@smooth (local or global) section}
 Let $M$ be a smooth manifold with or without boundary and let $\pi\colon \mathbf{E}\longrightarrow M$ be a smooth vector bundle. A \textbf{smooth (local or global) section of $\mathbf{E}$} is a (local or global) section that is a smooth map from its domain to $\mathbf{E}$.
 \end{definition}
 \begin{remark}\label{obs:nociones-fundamentales-seccion-local-o-global-en-bruto-de}
 For some purposes, we need to consider functions that satisfy the conditions for a section except for continuity. A function satisfying all the conditions for a section other than continuity will be called a \textit{\textbf{rough (local or global) section of $\mathbf{E}$}}.
 \end{remark}
\begin{figura}[h]
    \includegraphics[scale=1]{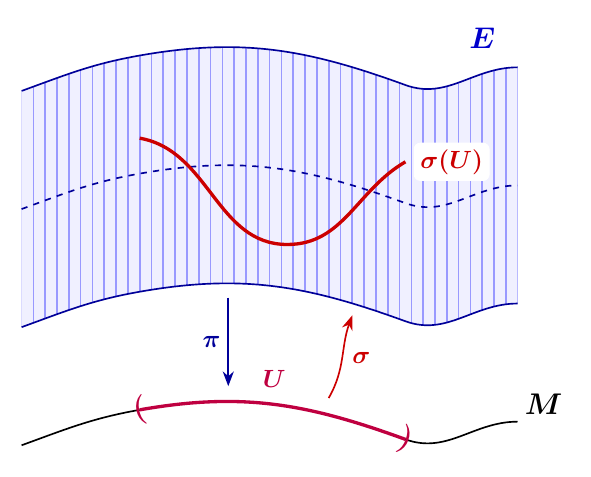}
    \caption{Local section of a vector bundle.}
    \label{fig:seccion-local-haz-vectorial}
\end{figura}
 Since we will consider sections of various degrees of regularity, it is useful to have notation that simplifies specifying the differentiability class of a section.
 \begin{definition}\label{secciones de clase Ck y de soporte compacto}\index{sections of class Ck and compactly supported sections@sections of class Ck and compactly supported sections}
 Let $M$ be a smooth manifold with or without boundary and let $\mathbf{E}\longrightarrow M$ be a smooth vector bundle. Let $k\in \mathbb{N}_{0}\cup \{\infty\}$. Define \[\Gamma^{k}(\mathbf{E}):=\{\boldsymbol{\sigma}\colon M\longrightarrow \mathbf{E}\mid \boldsymbol{\sigma}\text{ is a section of class $C^{k}$}\}.\]
 We will sometimes need to work with compactly supported sections, so we also define the support of a section. If $\boldsymbol{\sigma}\colon M\longrightarrow \mathbf{E}$ is a section (possibly rough), write
 \[\supp(\boldsymbol{\sigma}):=\overline{\{p\in M\mid \boldsymbol{\sigma}(p)\neq0_{p}\}}\] and also define \[\Gamma_{c}^{k}(\mathbf{E}):=\{\boldsymbol{\sigma}\in \Gamma^{k}(\mathbf{E})\mid \boldsymbol{\sigma}\text{ has compact support}\}.\] If $k=\infty$, we simply write \[\Gamma(\mathbf{E}):=\Gamma^{\infty}(\mathbf{E})\qquad \text{and}\qquad \Gamma_{c}(\mathbf{E}):=\Gamma_{c}^{\infty}(\mathbf{E}).\]
 When we need to specify the manifold on which sections are being considered, for example when working with several bundles over different smooth manifolds, we write $\Gamma^{k}(M,\mathbf{E})$, $\Gamma(M,\mathbf{E})$, $\Gamma_{c}^{k}(M,\mathbf{E})$, and $\Gamma_{c}(M,\mathbf{E})$ with the same meanings as above.
 \end{definition}
 We now define local and global frames. Among other uses, these will give an important criterion for determining when a section is smooth.
 \begin{definition}\label{def:nociones-fundamentales-marco-local-para-e-sobre-u}\index{local frame!of a vector bundle}
 Let $\pi\colon \mathbf{E}\longrightarrow M$ be a vector bundle, let $U\subseteq M$ be an open subset, and let $(\boldsymbol{\sigma}_{1},\dots,\boldsymbol{\sigma}_{k})$ be a $k$-tuple of local sections of $\mathbf{E}$ over $U$. We say that $(\boldsymbol{\sigma}_{1},\dots,\boldsymbol{\sigma}_{k})$ is a \textbf{local frame for E over U} if $\{\boldsymbol{\sigma}_{1}(p),\dots,\boldsymbol{\sigma}_{k}(p)\}$ is a basis for $\mathbf{E}_{p}$ for each $p\in U$. If $U=M$, the local frame $(\boldsymbol{\sigma}_{1},\dots,\boldsymbol{\sigma}_{k})$ is called a \textbf{global frame for $\mathbf{E}$ over $U$}. To simplify notation, we write $(\boldsymbol{\sigma}_{i})_{i=1}^{k}:=(\boldsymbol{\sigma}_{1},\dots,\boldsymbol{\sigma}_{k})$ or simply $(\boldsymbol{\sigma}_{i}):=(\boldsymbol{\sigma}_{1},\dots,\boldsymbol{\sigma}_{k})$ when the index set is clear.
 \end{definition}
 \begin{definition}\label{def:nociones-fundamentales-marco-local-suave}\index{smooth local frame}
 Let $M$ be a smooth manifold with or without boundary. If
 $\pi\colon \mathbf{E}\longrightarrow M$ is a smooth vector bundle,
 we say that a local frame $(\boldsymbol{\sigma}_{i})_{i=1}^{k}$ defined
 on an open subset $U$ of $M$ is a \textbf{smooth local frame} if
 $\boldsymbol{\sigma}_{i}$ is a smooth section for every
 $i\in \{1,\dots,k\}$.
 \end{definition}
 Smooth local frames and smooth local trivializations of a smooth vector bundle are related as follows:

 \begin{proposition}\label{marcos asociados con trivializaciones}
 Let $M$ be a smooth manifold with or without boundary. If
 $\pi\colon \mathbf{E}\longrightarrow M$ is a smooth vector bundle of
 rank $k$, then:
 \begin{enumerate}[label=(\alph*)]
 \item Every smooth local trivialization $\Phi\colon \pi^{-1}(U)\longrightarrow U\times\mathbb{R}^{k}$ induces a smooth local frame given by $\boldsymbol{\sigma}_{i}=\Phi^{-1}(p,\mathbf{e}_{i})$ for each $p\in U$.
 \item Every smooth local frame $(\boldsymbol{\sigma}_{i})$ defined on an open set $U\subseteq M$ induces a smooth local trivialization.
 \end{enumerate}
 \end{proposition}
 \begin{proof}
 \begin{enumerate}[label=(\alph*)]
 \item If $\Phi\colon \pi^{-1}(U)\longrightarrow U\times\mathbb{R}^{k}$ is a smooth local trivialization of $\mathbf{E}$, we can construct a smooth local frame for $\mathbf{E}$ over $U$. For each $i\in\{1,\dots,k\}$, define maps $\boldsymbol{\sigma}_{1},\dots,\boldsymbol{\sigma}_{k}\colon U\longrightarrow \mathbf{E}$ by $\boldsymbol{\sigma}_{i}(p)=\Phi^{-1}(p,\mathbf{e}_{i})=\Phi^{-1}\circ\widetilde{\mathbf{e}_{i}}(p)$, where $\widetilde{\mathbf{e}_{i}}\colon M\longrightarrow M\times\mathbb{R}^{k}$, $\widetilde{\mathbf{e}_{i}}(p)=(p,\mathbf{e}_{i})$, and $\{\mathbf{e}_i\}_{i=1}^{k}$ is the standard basis of $\mathbb{R}^{k}$.
 The family $(\widetilde{\mathbf{e}_{i}})_{i=1}^{k}$ forms a smooth global frame for the bundle $\pi_{1}\colon M\times\mathbb{R}^{k}\longrightarrow M$, where $\pi_{1}$ is projection onto the first coordinate.
 \[
 \begin{tikzcd}
 \pi^{-1}(U) \arrow{dr}{\pi} \arrow{rr}{\Phi}& & U\times \mathbb{R}^{k} \arrow{dl}[swap]{\pi_1} \\
 & \arrow[bend left = 30]{ul}{\boldsymbol{\sigma}_i} U \arrow[bend right = 30,swap]{ur}{\widetilde{\mathbf{e}}_i}
 \end{tikzcd}
 \]
 Then $\boldsymbol{\sigma}_{i}$ is smooth because $\Phi$ is a diffeomorphism; moreover, $\pi_{1}\circ\Phi=\pi$, so \[\pi\circ\boldsymbol{\sigma}_{i}(p)=\pi\circ \Phi^{-1}(p,\mathbf{e}_{i})=\pi_{1}(p,\mathbf{e}_{i})=p,\] which says that $\boldsymbol{\sigma}_{i}$ is a section. To see that $(\boldsymbol{\sigma}_{i}(p))$ forms a basis for $\mathbf{E}_{p}$, note that $\Phi$ restricts to an isomorphism from $\mathbf{E}_{p}$ to $\{p\}\times\mathbb{R}^{k}$ and that $\Phi(\boldsymbol{\sigma}_{i}(p))=(p,\mathbf{e}_{i})$; hence $\Phi$ sends $(\boldsymbol{\sigma}_{i}(p))$ to the standard basis for $\{p\}\times\mathbb{R}^{k}$. Therefore, $(\boldsymbol{\sigma}_{i})$ is a smooth local frame for $\mathbf{E}$, called \textbf{the smooth local frame associated with the smooth local trivialization $\Phi$}.
 \item Suppose $(\boldsymbol{\sigma}_{i})_{i=1}^{k}$ is a smooth local frame for $\mathbf{E}$ defined on an open set $U\subseteq M$. Define $\Psi\colon U\times\mathbb{R}^{k}\longrightarrow \pi^{-1}(U)$ by \[\Psi(p,(v^{1},\dots,v^{k}))=v^{i}\boldsymbol{\sigma}_{i}(p).\] Since $\{\boldsymbol{\sigma}_{i}(p)\mid i\in\{1,\dots,k\}\}$ forms a basis for $\mathbf{E}_{p}$ for each $p\in U$, $\Psi$ is bijective. Moreover, $\boldsymbol{\sigma}_{i}=\Psi\circ\widetilde{\mathbf{e}_{i}}$. If we prove that $\Psi$ is a diffeomorphism, then $\Psi^{-1}$ will be a smooth local trivialization whose associated smooth local frame is $(\boldsymbol{\sigma}_{i})^{k}_{i=1}$.
 Since $\Psi$ is bijective, it suffices to show that it is a local diffeomorphism. Let $q\in U$. By the definition of a smooth vector bundle, there exist an open set $V_{0}\subseteq M$ with $q\in V_{0}$ and a smooth local trivialization $\Phi\colon \pi^{-1}(V_{0})\longrightarrow V_{0}\times\mathbb{R}^{k}$. Replacing $V_{0}$ by $V:=V_{0}\cap U$ and restricting $\Phi$ gives a smooth local trivialization over an open set $V\subseteq U$ containing $q$. Since $\Phi$ is a diffeomorphism, it suffices to prove that
 \[
 (\Phi\circ\Psi)\restriction_{V\times\mathbb{R}^{k}}\colon V\times\mathbb{R}^{k}\longrightarrow V\times\mathbb{R}^{k}
 \]
 is a diffeomorphism.

 For each $i\in\{1,\dots,k\}$, the map $(\Phi\circ\boldsymbol{\sigma}_i)\restriction_{V}\colon V\longrightarrow V\times\mathbb{R}^{k}$ is smooth. Thus there exist smooth functions $\sigma_i^1,\dots,\sigma_i^k\colon V\longrightarrow\mathbb R$ such that $(\Phi\circ\boldsymbol{\sigma}_i)(p)=(p,(\sigma_i^1(p),\dots,\sigma_i^k(p)))$ for every $p\in V$. Consequently,
 \[
 (\Phi\circ\Psi)(p,(v^1,\dots,v^k))=(p,(v^i\sigma_i^1(p),\dots,v^i\sigma_i^k(p))),
 \]
 so $(\Phi\circ\Psi)\restriction_{V\times\mathbb{R}^{k}}$ is smooth.

 For each $p\in V$, the matrix $(\sigma_i^j(p))$ is invertible because $\{\boldsymbol{\sigma}_i(p)\mid i\in\{1,\dots,k\}\}$ is a basis for $\mathbf{E}_p$. Denote its inverse matrix by $(\tau_i^j(p))$. Since inversion on $GL(k,\mathbb R)$ is smooth, the functions $\tau_i^j$ are smooth and
 \[
 (\Phi\circ\Psi)^{-1}(p,(w^1,\dots,w^k))=(p,(w^i\tau_i^1(p),\dots,w^i\tau_i^k(p))).
 \]
 Therefore, $(\Phi\circ\Psi)^{-1}$ is smooth. This proves that $\Psi$ is a local diffeomorphism around $q$, and since $q\in U$ was arbitrary, the proof is complete.
 \end{enumerate}
 \end{proof}

\begin{remark}\label{cubierta marco local suave}
 In particular, Proposition~\ref{marcos asociados con trivializaciones} tells us that there always exists an open cover $(U_{\alpha})_{\alpha\in J}$ by domains of smooth local frames.
\end{remark}

 We obtain the following corollary:
 \begin{corollary}\label{carta del haz vectorial}
 Let $M$ be a smooth manifold with or without boundary, let
 $\pi\colon \mathbf{E}\longrightarrow M$ be a smooth vector bundle of rank
 $k$, let $(V,\phi)$ be a smooth chart on $M$ with coordinate functions
 $\phi=(x^{1},\dots,x^{n})$, and suppose we have a smooth local frame
 $(\boldsymbol{\sigma}_{i})_{i=1}^{k}$ for $\mathbf{E}$ over $V$.
 Define $\widetilde{\phi}\colon \pi^{-1}(V)\longrightarrow
 \phi(V)\times\mathbb{R}^{k}$ by
 \[\widetilde{\phi}(v^{i}\boldsymbol{\sigma}_{i})=(x^{1}(p),\dots,x^{n}(p),v^{1},\dots,v^{k}).\]
 Then $(\pi^{-1}(V),\widetilde{\phi})$ is a smooth coordinate chart of $\mathbf{E}$.
 \end{corollary}
\begin{proof}
The frame determines the smooth trivialization
\[
\Psi\colon V\times\mathbb R^k\longrightarrow\pi^{-1}(V),
\qquad
\Psi(p,v)=\displaystyle\sum_{i=1}^kv^i\boldsymbol\sigma_i(p),
\]
by Proposition~\ref{marcos asociados con trivializaciones}. The function
in the statement is
\(\widetilde\phi=(\phi\times\operatorname{Id}_{\mathbb R^k})\circ\Psi^{-1}\);
it is therefore a diffeomorphism onto
\(\phi(V)\times\mathbb R^k\).
\end{proof}
 We now have a smoothness criterion for rough local sections.
 \begin{proposition}[Regularity of a section in a local frame]
\label{suavidad de secciones en bruto locales}
Let $M$ be a smooth manifold with or without boundary,
$\pi\colon\mathbf E\to M$ a smooth vector bundle of rank $r$, and
$\boldsymbol\tau\colon U\to\mathbf E$ a raw section over an open set
$U\subseteq M$. If $(\mathbf e_a)_{a=1}^r$ is a smooth local frame on $U$ and
\[
 \boldsymbol\tau(p)=\sum_{a=1}^r\tau^a(p)\mathbf e_a(p),
\]
then, for $k\in\mathbb N_0\cup\{\infty\}$,
$\boldsymbol\tau$ is of class $C^k$ if and only if all the functions
$\tau^a$ are of class $C^k$. In particular, it is continuous if and only if
all its components are continuous.
For a global section, it suffices to check the condition on a cover
by domains of smooth frames.
\end{proposition}
\begin{proof}
By Proposition~\ref{marcos asociados con trivializaciones}, the frame
determines the diffeomorphism
\[
 \Psi\colon U\times\mathbb R^r\longrightarrow\pi^{-1}(U),
 \qquad \Psi(p,v)=\sum_{a=1}^r v^a\mathbf e_a(p).
\]
The representation of the section is
\[
 (\Psi^{-1}\circ\boldsymbol\tau)(p)
 =\bigl(p,(\tau^1(p),\dots,\tau^r(p))\bigr).
\]
Projection onto each component of the second entry proves necessity.
For sufficiency, this map into the product is continuous or of class
$C^k$ when every component is, and one composes with $\Psi$.
In base coordinates it is
$x\mapsto(x,\tau^1\circ\phi^{-1}(x),\dots,\tau^r\circ\phi^{-1}(x))$,
so the criterion includes boundary charts with the local-extension
convention. Finally, continuity and $C^k$ regularity are local
properties, so a cover by these domains suffices.
\end{proof}

The criterion is independent of the frame. Indeed, if another frame is
written
$\widetilde{\mathbf e}_b=\displaystyle\sum_{a=1}^r\mathbf e_a A^a{}_b$,
the matrix $A\colon U\to GL(r,\mathbb R)$ is smooth, and the component
columns satisfy
\begin{equation}\label{eq:compatibilidad-componentes-secciones}
 \widetilde\tau=A^{-1}\tau.
\end{equation}
The inverse is smooth because its entries are quotients of smooth
cofactors by the nonzero determinant. Thus multiplication by $A$ or
$A^{-1}$ preserves continuity and every class $C^k$.
When prescribing a section through component columns in different frames,
this relation must hold on overlaps; uniqueness of the coordinates of a
vector in a basis proves that the expressions glue to a unique raw section.

For a complex bundle of rank $r$, the same proof uses
$U\times\mathbb C^r$ in place of $U\times\mathbb R^r$.
A complex component function is continuous or of class $C^k$ when both
its real and imaginary parts are. The identification
$\mathbb C^r\cong\mathbb R^{2r}$ justifies this criterion, and the
change-of-frame matrices belong to $GL(r,\mathbb C)$.
This is regularity with respect to the real coordinates of the base.

\begin{example}[A scalar function as a section and a change of frame]
\label{ej:funciones-como-secciones-y-marco-movil}
In the trivial bundle $M\times\mathbb R\to M$, every raw section has the
form $\sigma_f(p)=(p,f(p))$. The frame $\mathbf e(p)=(p,1)$ shows that
$\sigma_f$ is continuous or $C^k$ exactly when $f$ is.
Thus the study of scalar functions is a special case of the study of
sections: the fiber is a copy of $\mathbb R$ at each point.

In $\mathbb R\times\mathbb R^2\to\mathbb R$, consider the smooth frame
\[
 \mathbf e_1(t)=(t,(\cos t,\sin t)),
 \qquad \mathbf e_2(t)=(t,(-\sin t,\cos t)).
\]
Its matrix relative to the constant frame is
\[
 A(t)=
 \begin{pmatrix}\cos t&-\sin t\\ \sin t&\cos t\end{pmatrix},
 \qquad
 A(t)^{-1}=
 \begin{pmatrix}\cos t&\sin t\\ -\sin t&\cos t\end{pmatrix}.
\]
The smooth section $\tau(t)=(t,(1,0))$ has components $(1,0)$ in the
constant frame and $(\cos t,-\sin t)$ in the moving frame:
\[
 (1,0)=\cos t\,(\cos t,\sin t)-\sin t\,(-\sin t,\cos t).
\]
The components change, while the vector in each fiber remains the same.
The regularity criterion gives the same answer in both frames.

The section $\eta(t)=(t,(|t|,0))$ is continuous and is not $C^1$ at $0$.
Replacing $|t|$ by $H(t):=\mathbf 1_{(0,\infty)}(t)$ gives a raw section
that is discontinuous at $0$. Projection onto the first fiber component
recovers $|t|$ or $H(t)$, so neither a smooth change of frame nor section
notation can remove the failure of regularity.
\end{example}
Multiplication by smooth functions makes the space of sections of a vector bundle a $C^{\infty}(M)$-module. This algebraic structure expresses the fact that the coefficients of a section can vary smoothly over the base.
\begin{proposition}\label{prop:nociones-fundamentales-haz-vectorial-suave-variedad-suave-frontera}
Let $\pi\colon \mathbf{E}\longrightarrow M$ be a smooth vector bundle over a smooth manifold $M$ with or without boundary. Define an action of $C^{\infty}(M)$ on $\Gamma(\mathbf{E})$ by
\[
(f\boldsymbol{\sigma})(p)\ :=\ f(p)\boldsymbol{\sigma}(p),\qquad f\in C^{\infty}(M),\ \boldsymbol{\sigma}\in\Gamma(\mathbf{E}),\ p\in M.
\]
Then $\Gamma(\mathbf{E})$, with the usual addition of functions and this action, is a $C^{\infty}(M)$-module.
\end{proposition}

\begin{proof}
Let $(U_{\alpha})_{\alpha\in J}$ be an open cover of $M$ such that, for each $\alpha\in J$, there exists a smooth local frame $(\boldsymbol{\sigma}_{\alpha,i})_{i=1}^{k}$ of $\mathbf{E}$ defined on $U_{\alpha}$.
If $\boldsymbol{\sigma}\in\Gamma(\mathbf{E})$, on each $U_{\alpha}$ we can write
\[
\boldsymbol{\sigma}\restriction_{U_{\alpha}}
=\displaystyle\sum_{i=1}^{k}a_{\alpha}^{i}\boldsymbol{\sigma}_{\alpha,i},
\qquad a_{\alpha}^{i}\in C^{\infty}(U_{\alpha}),
\]
where $k$ is the rank of $\mathbf{E}$. For $f\in C^{\infty}(M)$, we have
\[
(f\boldsymbol{\sigma})\restriction_{U_{\alpha}}
=\displaystyle\sum_{i=1}^{k}\bigl(f\restriction_{U_{\alpha}}\cdot a_{\alpha}^{i}\bigr)\boldsymbol{\sigma}_{\alpha,i}.
\]
By Proposition~\ref{suavidad de secciones en bruto locales}, $f\boldsymbol{\sigma}$ is smooth on $U_{\alpha}$ because $f\restriction_{U_{\alpha}}\cdot a_{\alpha}^{i}$ is smooth for every $i$. Since this holds for every $\alpha\in J$, we conclude that $f\boldsymbol{\sigma}\in\Gamma(\mathbf{E})$.

The module identities are
\[
(f+g)\boldsymbol{\sigma}=f\boldsymbol{\sigma}+g\boldsymbol{\sigma},\qquad
f(\boldsymbol{\sigma}+\boldsymbol{\tau})=f\boldsymbol{\sigma}+f\boldsymbol{\tau},\qquad
(fg)\boldsymbol{\sigma}=f(g\boldsymbol{\sigma}),\qquad
1\cdot\boldsymbol{\sigma}=\boldsymbol{\sigma},
\]
To verify them, evaluate each side at a point $p\in M$ and use
the vector space identities in the fiber $\mathbf E_p$. For
example,
\[
 ((fg)\boldsymbol{\sigma})(p)
 =f(p)g(p)\boldsymbol{\sigma}(p)
 =\bigl(f(g\boldsymbol{\sigma})\bigr)(p).
\]
The other three equalities follow in the same way. Therefore,
$\Gamma(\mathbf{E})$ is a $C^{\infty}(M)$\nobreakdash-module.
\end{proof}
The natural notion of a map between two vector bundles should cover a
map between their base spaces and be linear on each fiber. This leads to
the following definition.
\begin{definition}\label{def:nociones-fundamentales-homomorfismo-de-haces}\index{bundle homomorphism}
 If $\pi\colon \mathbf{E}\longrightarrow M$ and $\pi'\colon \mathbf{E}'\longrightarrow M'$ are vector bundles over the same field $\mathbb K\in\{\mathbb R,\mathbb C\}$, a continuous function $F\colon \mathbf{E}\longrightarrow \mathbf{E}'$ is called a \textit{\textbf{bundle homomorphism}} if there exists a function $f\colon M\longrightarrow M'$ such that $\pi'\circ F=f\circ \pi$, that is, such that the following diagram commutes

\begin{center}
\begin{tikzcd}
\mathbf{E} \arrow[r, "F"] \arrow[d, "\pi"] & \mathbf{E}' \arrow[d, "\pi'"] \\
M \arrow[r, "f"'] & M'
\end{tikzcd}

\end{center}
with the property that, for every $p\in M$, $F\restriction_{\mathbf{E}_{p}}\colon \mathbf{E}_{p}\longrightarrow \mathbf{E}'_{f(p)}$ is $\mathbb K$-linear. We express the relation between $F$ and $f$ by saying that \textit{\textbf{$F$ covers $f$}}.

If $\mathbf{E}$, $\mathbf{E}'$, and $F$ are smooth, we say that $F$ is a
\textit{\textbf{smooth bundle homomorphism}}. When both bundles have
the same base and the covered map is the identity, we write
\[
\underline{\operatorname{Hom}}(\mathbf{E},\mathbf{E}')
:=\{F\colon \mathbf{E}\longrightarrow \mathbf{E}'\mid
F\text{ is a smooth bundle homomorphism over }M\}.
\]
\end{definition}
We have not required $f$ to be continuous; this follows from continuity of $F$, which also determines it uniquely.

\begin{proposition}\label{prop:nociones-fundamentales-haces-vectoriales-homomorfismo-haces-cubre-continua}
 Let $\pi\colon \mathbf{E}\longrightarrow M$ and
 $\pi'\colon \mathbf{E}'\longrightarrow M'$ be vector bundles, and let
 $F\colon \mathbf{E}\longrightarrow \mathbf{E}'$ be a bundle
 homomorphism covering $f\colon M\longrightarrow M'$. Then $f$ is
 continuous and uniquely determined by $F$. If, in addition, $M$ and $M'$
 are smooth manifolds with or without boundary, the bundles are smooth, and $F$ is
 smooth, then $f$ is also smooth.
\end{proposition}
\begin{proof}
 We have $f=\pi'\circ F\circ \zeta$, where $\zeta\colon M\longrightarrow \mathbf{E}$ is such that $\zeta(p)=0_{\mathbf{E}_{p}}\in \mathbf{E}_{p}$ for each $p\in M$. To verify this, note that $(\pi'\circ F\circ \zeta)(p)=\pi'(F(0_{\mathbf{E}_{p}}))=\pi'(0_{\mathbf{E}_{f(p)}})=f(p)$. Thus $f$ is continuous, and if the bundles and $F$ are smooth, then $\pi'$ and $\zeta$ are smooth, so $f$ is smooth.

 Moreover, $f$ is uniquely determined by $F$, since every $f$ making the diagram commute must equal that composition.
\end{proof}

 \begin{definition}\label{def:nociones-fundamentales-homomorfismo-haces-biyectivo-cuya-inversa-homomorfismo}\index{bundle isomorphism}
 If $F\colon \mathbf{E}\longrightarrow \mathbf{E}'$ is a bijective bundle homomorphism whose inverse is a bundle homomorphism, we say that $F$ is a \textit{bundle isomorphism}. If $F$ is a diffeomorphism, we say that $F$ is a \textit{smooth bundle isomorphism}. \end{definition}

The next exercise shows that, over the same base,
bijectivity of a smooth homomorphism suffices for its inverse to be smooth.
\begin{exercise}\label{homomorfismo biyectivo es isomorfismo}
 Let $M$ be a smooth manifold with or without boundary. Prove that if $\mathbf{E}$ and $\mathbf{E}'$ are smooth vector bundles over $M$ and $F\in \underline{\operatorname{Hom}}(\mathbf{E},\mathbf{E}')$ is bijective, then it is a smooth bundle isomorphism.
\end{exercise}
 \begin{remark}\label{obs F homomorfismo es C infty lineal}
 If $\mathbf{E},\mathbf{E}'\longrightarrow M$ are smooth vector bundles over a smooth manifold $M$ with or without boundary and $F\in \underline{\operatorname{Hom}}(\mathbf{E},\mathbf{E}')$, then composition with $F$ induces a map \[\widetilde{F}\colon \Gamma(\mathbf{E})\longrightarrow \Gamma(\mathbf{E}')\] as follows:

 \[\widetilde{F}(\boldsymbol{\sigma})(p)=(F\circ \boldsymbol{\sigma})(p)=F(\boldsymbol{\sigma}(p)).\]
 We have $\widetilde{F}(\boldsymbol{\sigma})\in \Gamma(\mathbf{E}')$, since it is smooth as a composition of smooth functions. If the bundles are defined over $\mathbb K\in\{\mathbb R,\mathbb C\}$, linearity of $F$ on each fiber implies that $\widetilde F$ is $\mathbb K$-linear on sections. Moreover, for $f\in C^\infty(M)$,
 \[
 \widetilde F(f\boldsymbol\sigma)(p)
 =F_p\bigl(f(p)\boldsymbol\sigma(p)\bigr)
 =f(p)F_p\bigl(\boldsymbol\sigma(p)\bigr)
 =\bigl(f\widetilde F(\boldsymbol\sigma)\bigr)(p),
 \]
 so $\widetilde{F}$ is $C^{\infty}(M)$-linear; that is, $\widetilde{F}\in \displaystyle\operatorname{Hom}_{C^{\infty}(M)}(\Gamma(\mathbf{E}),\Gamma(\mathbf{E}'))$.

 \end{remark}

The preceding remark tells us that every bundle homomorphism induces a homomorphism between the modules of smooth sections. The converse will also hold: every module homomorphism between the sections of the bundles induces a bundle homomorphism.
To prove this, we first need a lemma.
\begin{lemma}[Section with prescribed value in a fiber]\label{lema:seccion-valor-prescrito}\index{section with prescribed value in a fiber@section with prescribed value in a fiber}
Let $M$ be a smooth manifold with or without boundary, let
$\pi\colon \mathbf{E}\longrightarrow M$ be a smooth vector bundle, let
$p\in M$, and let $v\in \mathbf{E}_p$. Then there exists a smooth global
section $\boldsymbol{\sigma}\in\Gamma(\mathbf{E})$ with
$\boldsymbol{\sigma}(p)=v$. Moreover, if $U$ is a neighborhood of $p$ over
which $\mathbf{E}$ admits a local trivialization, we can choose
$\boldsymbol{\sigma}$ with $\operatorname{supp}\boldsymbol{\sigma}\subset U$.
\end{lemma}

\begin{proof}
Let $(U,\Phi)$ be a local trivialization of $\mathbf{E}$ on an open neighborhood of $p$, with
\[
\Phi\colon \pi^{-1}(U)\longrightarrow U\times\mathbb{R}^k.
\]
Write $\Phi(v)=(p,u_0)$ with $u_0\in\mathbb{R}^k$. By Proposition~\ref{funciones flan}, we can choose a bump function $\chi\in C^\infty(M)$ with $\chi(p)=1$ and $\operatorname{supp}\chi\subset U$. Define
\[
\boldsymbol{\sigma}(q)=
\begin{cases}
\Phi^{-1}\bigl(q,\chi(q)u_0\bigr), & q\in U,\\[2pt]
0_q, & q\notin U,
\end{cases}
\]
where $0_q$ denotes the zero vector in the fiber $\mathbf{E}_q$. On $U$, the composition $q\mapsto \chi(q)u_0$ is smooth and, since $\Phi^{-1}$ is smooth, $\boldsymbol{\sigma}\restriction_{U}$ is smooth. If $q\notin \overline{U}$, then $\boldsymbol{\sigma}$ agrees with the zero section, which is smooth.

To justify gluing at the boundary, we use the condition $\operatorname{supp}\chi\subset U$. By definition, $M\setminus\operatorname{supp}\chi$ is open and contains $\partial U$, so for each $q_0\in\partial U$ there exists an open neighborhood $V_{q_0}\subset M$ such that $V_{q_0}\cap\operatorname{supp}\chi=\varnothing$, and in particular $\chi\restriction_{V_{q_0}}=0$. Consequently, for $q\in V_{q_0}\cap U$ we have
\[
\boldsymbol{\sigma}(q)=\Phi^{-1}\bigl(q,\chi(q)u_0\bigr)=\Phi^{-1}(q,0_{\mathbb{R}^k})=0_q,
\]
and for $q\in V_{q_0}\setminus U$ we also have $\boldsymbol{\sigma}(q)=0_q$ by definition. In other words, on an open neighborhood of each point of $\partial U$, both expressions for $\boldsymbol{\sigma}$ agree exactly with the zero section. Therefore, $\boldsymbol{\sigma}$ is smooth on $M$.

By construction, $\boldsymbol{\sigma}(p)=\Phi^{-1}(p,u_0)=v$. Moreover, if $q\in U$ and $\chi(q)=0$, then $\boldsymbol{\sigma}(q)=0_q$, and if $q\notin U$, we also have $\boldsymbol{\sigma}(q)=0_q$, so $\operatorname{supp}\boldsymbol{\sigma}\subset\operatorname{supp}\chi\subset U$.
\end{proof}

\begin{lemma}[Characterization of bundle homomorphisms]\label{lema de caracterizacion de homomorfismos de haces}\index{characterization of bundle homomorphisms@characterization of bundle homomorphisms}
Let $M$ be a smooth manifold with or without boundary, and let
$\pi\colon \mathbf{E}\longrightarrow M$ and
$\pi'\colon \mathbf{E}'\longrightarrow M$ be smooth vector bundles. A map
\[
\mathcal{F}\colon \Gamma(\mathbf{E})\longrightarrow \Gamma(\mathbf{E}')
\]
is $C^{\infty}(M)$-linear if and only if there exists a \emph{unique} smooth bundle homomorphism $F\colon \mathbf{E}\longrightarrow \mathbf{E}'$ over $M$ such that
\[
\mathcal{F}(\boldsymbol{\sigma})=F\circ \boldsymbol{\sigma},\qquad \forall\boldsymbol{\sigma}\in\Gamma(\mathbf{E}).
\]
\end{lemma}
\begin{proof}
If $F\colon \mathbf{E}\to \mathbf{E}'$ is a smooth homomorphism over $M$, then
$\boldsymbol{\sigma}\mapsto F\circ\boldsymbol{\sigma}$ is $C^\infty(M)$-linear by linearity
of each $F_p$.

Conversely, let
$\mathcal F\colon\Gamma(\mathbf{E})\to\Gamma(\mathbf{E}')$ be a
$C^\infty(M)$-linear map. For $p\in M$ and $v\in \mathbf{E}_p$, choose
$\boldsymbol{\sigma}\in\Gamma(\mathbf{E})$ with $\boldsymbol{\sigma}(p)=v$, whose existence is guaranteed
by Lemma~\ref{lema:seccion-valor-prescrito}, and define
\[
 F_p(v):=\mathcal F(\boldsymbol{\sigma})(p).
\]
We show that the value is independent of $\boldsymbol{\sigma}$. It suffices to consider a section
$\boldsymbol{\tau}$ with $\boldsymbol{\tau}(p)=0$ and verify that
$\mathcal F(\boldsymbol{\tau})(p)=0$. Take a
local frame $(\mathbf{e}_1,\ldots,\mathbf{e}_r)$ of $\mathbf{E}$ on an open set $U$ containing
$p$, and a bump function $\chi$ supported in $U$ and equal to
one in a neighborhood of $p$. On $U$, write
$\displaystyle \boldsymbol{\tau}=\displaystyle\sum_{i=1}^{r} a^i \mathbf{e}_i$, where $a^i(p)=0$. Also choose a bump
function $\rho$ supported in $U$ and equal to one on
$\operatorname{supp}\chi$. Extend
$\widetilde a^i:=\chi a^i$ and $\widetilde{\mathbf{e}}_i:=\rho \mathbf{e}_i$ by zero. These are,
respectively, smooth global functions and sections,
$\widetilde a^i(p)=0$, and
$\displaystyle \chi\boldsymbol{\tau}=\displaystyle\sum_{i=1}^{r}\widetilde a^i\widetilde{\mathbf{e}}_i$.
$C^\infty(M)$-linearity gives
\[
 \mathcal F(\chi\boldsymbol{\tau})(p)
 =\displaystyle\sum_{i=1}^{r}\widetilde a^i(p)\mathcal F(\widetilde{\mathbf{e}}_i)(p)=0.
\]
The section $(1-\chi)\boldsymbol{\tau}$ vanishes in a neighborhood of $p$. If
$\eta$ is a bump function with $\eta(p)=1$ and support contained in that
neighborhood, then
\[
 \mathcal F((1-\chi)\boldsymbol{\tau})(p)
 =\eta(p)\mathcal F((1-\chi)\boldsymbol{\tau})(p)
 =\mathcal F(\eta(1-\chi)\boldsymbol{\tau})(p)=0.
\]
Thus $\mathcal F(\boldsymbol{\tau})(p)=0$, and $F_p$ is well defined. Linearity of
$F_p$ follows by choosing sections with prescribed values and using
linearity of $\mathcal F$.

Define $F\colon \mathbf{E}\to \mathbf{E}'$ by
$F(v)=F_{\pi(v)}(v)$. To prove smoothness, let
$(\mathbf{e}_1,\ldots,\mathbf{e}_r)$ and
$(\mathbf{e}'_1,\ldots,\mathbf{e}'_s)$ be local frames defined on an
open set $U$ around $p$. Choose an open set $V$ such that
$p\in V$ and $\overline V\subseteq U$, and a bump function
$\rho\in C^\infty(M)$ with $\rho=1$ in a neighborhood of $\overline V$ and
$\operatorname{supp}\rho\subseteq U$. For each $i$, define
\[
 \widetilde{\mathbf e}_i(q)=
 \begin{cases}
  \rho(q)\mathbf e_i(q),&q\in U,\\[2pt]
  0_q,&q\notin U.
 \end{cases}
\]
The support condition on $\rho$ shows, as in
Lemma~\ref{lema:seccion-valor-prescrito}, that
$\widetilde{\mathbf e}_i$ is a smooth global section; moreover, it agrees
with $\mathbf e_i$ on $V$. There exist smooth functions $A^a_i$ such that
\[
 \mathcal F(\widetilde{\mathbf{e}}_i)=\displaystyle\sum_{a=1}^s A^a_i \mathbf{e}'_a
\]
on $V$. If $q\in V$ and $\displaystyle v=\displaystyle\sum_{i=1}^{r} z^i \mathbf{e}_i(q)$, the independence just
proved gives
\[
 F_q(v)=\displaystyle\sum_{a=1}^s\displaystyle\sum_{i=1}^r A^a_i(q)z^i \mathbf{e}'_a(q).
\]
In trivializations, the total map thus has the expression
$(q,z)\mapsto(q,A(q)z)$, which is smooth. By construction, it covers the identity,
is linear on each fiber, and satisfies
$F\circ\boldsymbol{\sigma}=\mathcal F(\boldsymbol{\sigma})$ for every $\boldsymbol{\sigma}\in\Gamma(\mathbf{E})$.

If $G$ is another homomorphism with the same action on sections, given
$p\in M$ and $v\in \mathbf{E}_p$, choose $\boldsymbol{\sigma}(p)=v$ and obtain
$F_p(v)=(F\circ\boldsymbol{\sigma})(p)=(G\circ\boldsymbol{\sigma})(p)=G_p(v)$. Therefore, $F=G$.
\end{proof}
The spaces of linear maps between the fibers themselves form a
vector bundle. The following proposition makes this construction
precise.

\begin{proposition}\label{haz Hom vectorial}
Let $M$ be a smooth manifold with or without boundary, and let
$\mathbf{E}\to M$ and $\mathbf{E}'\to M$ be smooth vector bundles of ranks
$k$ and $l$. Define
\[
\operatorname{Hom}(\mathbf{E},\mathbf{E}') \ :=\ \coprod_{p\in M}\operatorname{Hom}(\mathbf{E}_p,\mathbf{E}'_p),
\qquad
\pi_{\operatorname{Hom}}\colon \operatorname{Hom}(\mathbf{E},\mathbf{E}')\longrightarrow M,\quad (p,T_p)\mapsto p,
\]
with $T_p\in \operatorname{Hom}(\mathbf{E}_p,\mathbf{E}'_p)$. Then $\operatorname{Hom}(\mathbf{E},\mathbf{E}')$ admits a unique topology and a unique smooth structure making it a smooth vector bundle of rank $kl$ over $M$ for which the matrix trivializations induced by smooth local frames of $\mathbf E$ and $\mathbf E'$ are smooth.
\end{proposition}

\begin{proof}
Let $(U_{\alpha})_{\alpha\in J}$ be an open cover of $M$ such that, for every $\alpha\in J$, there exist smooth local frames
\[
(\boldsymbol{\sigma}_{\alpha,i})_{i=1}^{k}\ \text{of }\mathbf{E}\ \text{on }U_{\alpha},
\qquad
(\boldsymbol{\sigma}'_{\alpha,j})_{j=1}^{l}\ \text{of }\mathbf{E}'\ \text{on }U_{\alpha}.
\]
The associated trivializations
\[
\Phi^{\mathbf{E}}_{\alpha}\colon \pi_{\mathbf{E}}^{-1}(U_{\alpha})\longrightarrow U_{\alpha}\times \mathbb{R}^{k},
\qquad
\Phi^{\mathbf{E}'}_{\alpha}\colon \pi_{\mathbf{E}'}^{-1}(U_{\alpha})\longrightarrow U_{\alpha}\times \mathbb{R}^{l}
\]
are characterized by
\[
\Phi^{\mathbf{E}}_{\alpha}\bigl(\boldsymbol{\sigma}_{\alpha,i}(p)\bigr)=(p,e_i),\qquad
\Phi^{\mathbf{E}'}_{\alpha}\bigl(\boldsymbol{\sigma}'_{\alpha,j}(p)\bigr)=(p,e'_j),
\]
where $(e_i)_{i=1}^{k}$ and $(e'_j)_{j=1}^{l}$ are the standard bases of $\mathbb{R}^{k}$ and $\mathbb{R}^{l}$.

Given $p\in U_{\alpha}$ and $T_p\in \operatorname{Hom}(\mathbf{E}_p,\mathbf{E}'_p)$, write its matrix in the local frames as
\[
T_{\alpha}(p)=\bigl(T_{\alpha}(p)^{j}_{i}\bigr)_{j,i}\in M_{l\times k}(\mathbb{R}),
\]
defined by
\[
T_p(\boldsymbol{\sigma}_{\alpha,i}(p))=\displaystyle\sum_{j=1}^{l}T_{\alpha}(p)^{j}_{i}\boldsymbol{\sigma}'_{\alpha,j}(p).
\]
In particular, the linear map \(T_{\alpha}(p)\colon \mathbb{R}^{k}\longrightarrow\mathbb{R}^{l}\) satisfies
\[
T_{\alpha}(p)(e_i)=\displaystyle\sum_{j=1}^{l}T_{\alpha}(p)^{j}_{i}e'_j.
\]

Then define the local trivialization of $\operatorname{Hom}(\mathbf{E},\mathbf{E}')$ by
\[
\Phi^{\operatorname{Hom}}_{\alpha}(T_p):=(p,T_{\alpha}(p)).
\]
This function is bijective, with inverse given by
\[
(p,S)\ \longmapsto\ \left(u=\displaystyle\sum_{i=1}^{k}u^{i}\boldsymbol{\sigma}_{\alpha,i}(p)\ \mapsto\
\displaystyle\sum_{j=1}^{l}\displaystyle\sum_{i=1}^{k}S^{j}_{i}u^{i}\boldsymbol{\sigma}'_{\alpha,j}(p)\right),
\qquad S\in M_{l\times k}(\mathbb{R}).
\]

Here $S$ is an arbitrary matrix in $M_{l\times k}(\mathbb{R})$.
If $S=T_\alpha(p)$, the reconstructed map agrees exactly with $T_p$;
and if we start from an arbitrary $S$ and then apply $\Phi^{\operatorname{Hom}}_{\alpha}$, we recover $(p,S)$.
This shows that the two functions are indeed inverse to each other.

In particular, for each $p\in U_{\alpha}$, the restriction
\[
\Phi^{\operatorname{Hom}}_{\alpha}\restriction_{\operatorname{Hom}(\mathbf{E}_p,\mathbf{E}'_p)}\colon \operatorname{Hom}(\mathbf{E}_p,\mathbf{E}'_p)\longrightarrow \{p\}\times M_{l\times k}(\mathbb{R})
\]
is a linear isomorphism of real vector spaces.

If $U_{\alpha}\cap U_{\beta}\neq\varnothing$, Lemma~\ref{lema existnecia mapeo GL} says that the changes of trivialization of $\mathbf{E}$ and $\mathbf{E}'$ are described by smooth functions
\[
\tau^{\mathbf{E}}_{\alpha\beta}\colon U_{\alpha}\cap U_{\beta}\longrightarrow GL(k,\mathbb{R}),\qquad
\tau^{\mathbf{E}'}_{\alpha\beta}\colon U_{\alpha}\cap U_{\beta}\longrightarrow GL(l,\mathbb{R}),
\]
which, by Proposition~\ref{marcos asociados con trivializaciones}, are characterized by
\[
\boldsymbol{\sigma}_{\beta,i}=\displaystyle\sum_{m=1}^{k}(\tau^{\mathbf{E}}_{\alpha\beta})^{m}_{i}\boldsymbol{\sigma}_{\alpha,m},\qquad
\boldsymbol{\sigma}'_{\beta,j}=\displaystyle\sum_{n=1}^{l}(\tau^{\mathbf{E}'}_{\alpha\beta})^{n}_{j}\boldsymbol{\sigma}'_{\alpha,n}.
\]

With this notation, if the matrix of $T_p$ in $U_{\beta}$ is $T_{\beta}(p)$, its matrix in $U_{\alpha}$ is
\[
T_{\alpha}(p)=\tau^{\mathbf{E}'}_{\alpha\beta}(p)T_{\beta}(p)\bigl(\tau^{\mathbf{E}}_{\alpha\beta}(p)\bigr)^{-1}.
\]

Consequently, the change of trivialization of $\operatorname{Hom}(\mathbf{E},\mathbf{E}')$ is described by the function
\[
\tau^{\operatorname{Hom}}_{\alpha\beta}\colon U_{\alpha}\cap U_{\beta}\longrightarrow GL(kl,\mathbb{R}),
\]
defined at each $p\in U_{\alpha}\cap U_{\beta}$ by
\[
\tau^{\operatorname{Hom}}_{\alpha\beta}(p)(S)=\tau^{\mathbf{E}'}_{\alpha\beta}(p)S(\tau^{\mathbf{E}}_{\alpha\beta}(p))^{-1},\qquad S\in M_{l\times k}(\mathbb{R}).
\]
Indeed,
\[(\Phi_{\alpha}^{\operatorname{Hom}}\circ(\Phi_{\beta}^{\operatorname{Hom}})^{-1})(p,S)=\Phi_{\alpha}^{\operatorname{Hom}}\left(u=\displaystyle\sum_{i=1}^{k}u^{i}\boldsymbol{\sigma}_{\beta,i}(p)\ \mapsto\
\displaystyle\sum_{j=1}^{l}\displaystyle\sum_{i=1}^{k}S^{j}_{i}u^{i}\boldsymbol{\sigma}'_{\beta,j}(p)\right)=(p,T_{\alpha}(p))\]
\[=(p,\tau^{\mathbf{E}'}_{\alpha\beta}(p)S\left(\tau^{\mathbf{E}}_{\alpha\beta}(p)\right)^{-1}).\]
For each $p$, the map
$S\mapsto\tau^{\operatorname{Hom}}_{\alpha\beta}(p)(S)$ is linear. Moreover,
it depends smoothly on $p$, since matrix multiplication and inversion
are smooth and the functions $\tau^{\mathbf{E}}_{\alpha\beta}$ and
$\tau^{\mathbf{E}'}_{\alpha\beta}$ are smooth.

Consequently, the transition functions of
$\operatorname{Hom}(\mathbf{E},\mathbf{E}')$ are smooth and linear on each
fiber. Lemma~\ref{lema del haz vectorial suave} provides the unique
topology and smooth structure making this space a smooth vector
bundle of rank $kl$ over $M$ for which the maps
$\Phi_\alpha^{\operatorname{Hom}}$ are smooth local trivializations.
The same transition calculation applies to any other pair of
smooth local frames of the original bundles. Their matrix
trivializations are therefore smooth for this structure, and the asserted uniqueness
follows from uniqueness for the family already constructed.
\end{proof}

We use the notation $\operatorname{End}(\mathbf{E}):=\operatorname{Hom}(\mathbf{E},\mathbf{E})$.

The concepts introduced above are related as follows:

\begin{proposition}\label{iso Hom-Gamma-Hom}
Let $M$ be a smooth manifold with or without boundary, and let $\pi\colon \mathbf{E}\longrightarrow M$ and $\pi'\colon \mathbf{E}'\longrightarrow M$ be smooth vector bundles. Then we have the following natural identifications:

\begin{enumerate}[label=(\alph*)]
\item Each smooth bundle homomorphism $F\colon \mathbf{E}\longrightarrow \mathbf{E}'$ over $M$ determines a smooth section $p\longmapsto \bigl(F_p\colon \mathbf{E}_p\longrightarrow \mathbf{E}'_p\bigr)$ of the bundle $\operatorname{Hom}(\mathbf{E},\mathbf{E}')$. Conversely, every section $\mathbf{A}\in\Gamma(\operatorname{Hom}(\mathbf{E},\mathbf{E}'))$ defines a bundle homomorphism $F\colon \mathbf{E}\longrightarrow \mathbf{E}'$ by $F(v)=\mathbf{A}(\pi(v))(v)$ for every $v\in \mathbf{E}$. This gives a natural bijection
\[
\underline{\operatorname{Hom}}(\mathbf{E},\mathbf{E}')\cong \Gamma\bigl(\operatorname{Hom}(\mathbf{E},\mathbf{E}')\bigr),
\]
where $\underline{\operatorname{Hom}}(\mathbf{E},\mathbf{E}')$ denotes the set of bundle homomorphisms over $M$.

\item Every section $\mathbf{A}\in\Gamma(\operatorname{Hom}(\mathbf{E},\mathbf{E}'))$ induces a map $\mathcal{F}_{\mathbf{A}}\colon \Gamma(\mathbf{E})\longrightarrow\Gamma(\mathbf{E}')$ given by $\bigl(\mathcal{F}_{\mathbf{A}}(\boldsymbol{\sigma})\bigr)(p)=\mathbf{A}(p)(\boldsymbol{\sigma}(p))$ for $\boldsymbol{\sigma}\in\Gamma(\mathbf{E})$. This map is $C^{\infty}(M)$-linear. Conversely, if $\mathcal{F}\colon \Gamma(\mathbf{E})\longrightarrow\Gamma(\mathbf{E}')$ is $C^{\infty}(M)$-linear, there exists a unique section $\mathbf{A}\in\Gamma(\operatorname{Hom}(\mathbf{E},\mathbf{E}'))$ such that $\mathcal{F}=\mathcal{F}_{\mathbf{A}}$. This gives a natural isomorphism of $C^{\infty}(M)$-modules
\[
\Gamma\bigl(\operatorname{Hom}(\mathbf{E},\mathbf{E}')\bigr)\cong \operatorname{Hom}_{C^{\infty}(M)}\bigl(\Gamma(\mathbf{E}),\Gamma(\mathbf{E}')\bigr).
\]

\item Composing the two correspondences above gives a natural identification
\[
\underline{\operatorname{Hom}}(\mathbf{E},\mathbf{E}')\cong \operatorname{Hom}_{C^{\infty}(M)}\bigl(\Gamma(\mathbf{E}),\Gamma(\mathbf{E}')\bigr),
\]
compatible with composition of bundle homomorphisms and composition of $C^{\infty}(M)$-linear maps.
\end{enumerate}
\end{proposition}

\begin{proof}
(a) Let $F\colon \mathbf{E}\longrightarrow \mathbf{E}'$ be a bundle homomorphism over $M$, that is, a smooth map such that $\pi'\circ F=\pi$. For each $p\in M$, the restriction $F_p\colon \mathbf{E}_p\longrightarrow \mathbf{E}'_p$ is linear. Define a section $\mathbf{A}_F\in\Gamma(\operatorname{Hom}(\mathbf{E},\mathbf{E}'))$ by $\mathbf{A}_F(p):=F_p$. To see that $\mathbf{A}_F$ is smooth, it suffices to work in local trivializations.

On an open set $U_{\alpha}$ with smooth local frames $(\boldsymbol{\sigma}_{\alpha,i})_{i=1}^{k}$ of $\mathbf{E}$ and $(\boldsymbol{\sigma}'_{\alpha,j})_{j=1}^{l}$ of $\mathbf{E}'$, every vector $v\in \mathbf{E}_p$ can be written as $v=\displaystyle\sum\limits_{i=1}^{k}u^i\boldsymbol{\sigma}_{\alpha,i}(p)$. Denote by $T_{\alpha}(p)=\bigl(T_{\alpha}(p)^{j}_{i}\bigr)\in M_{l\times k}(\mathbb{R})$ the \emph{matrix} of $F_p$ with respect to these frames; that is,
\[
F_p\bigl(\boldsymbol{\sigma}_{\alpha,i}(p)\bigr)=\displaystyle\sum\limits_{j=1}^{l}T_{\alpha}(p)^{j}_{i}\boldsymbol{\sigma}'_{\alpha,j}(p)
\quad\text{for each }i\in\{1,\dots,k\}.
\]
Then, for $v\in \mathbf{E}_p$,
\[
F_p(v)=\displaystyle\sum\limits_{i=1}^{k}u^iF_p\bigl(\boldsymbol{\sigma}_{\alpha,i}(p)\bigr)
 =\displaystyle\sum\limits_{i=1}^{k}u^i\displaystyle\sum\limits_{j=1}^{l}T_{\alpha}(p)^{j}_{i}\boldsymbol{\sigma}'_{\alpha,j}(p).
\]

Since $F$ is smooth, the associated matrix $T_{\alpha}(p)$ depends smoothly on $p$, and by Proposition~\ref{suavidad de secciones en bruto locales}, $\mathbf{A}_F$ is a smooth section of $\operatorname{Hom}(\mathbf{E},\mathbf{E}')$.

Conversely, let $\mathbf{A}\in\Gamma(\operatorname{Hom}(\mathbf{E},\mathbf{E}'))$. For each $p\in M$ and $v\in \mathbf{E}_p$, define $F_{\mathbf{A}}(v):=\mathbf{A}(p)(v)$. Then $F_{\mathbf{A}}\colon \mathbf{E}\longrightarrow \mathbf{E}'$ is a smooth map, linear on each fiber, and moreover $\pi'(F_{\mathbf{A}}(v))=p=\pi(v)$, so $F_{\mathbf{A}}$ is a smooth bundle homomorphism over $M$.

If we start with $F$ and construct $\mathbf{A}_F$, then $F_{\mathbf{A}_F}(v)=\mathbf{A}_F(p)(v)=F_p(v)=F(v)$ for every $v\in \mathbf{E}_p$, so $F_{\mathbf{A}_F}=F$.

In the other direction, if we start with $\mathbf{A}$ and construct $F_{\mathbf{A}}$, for each
$p\in M$ we have
$\mathbf{A}_{F_{\mathbf{A}}}(p)=(F_{\mathbf{A}})_p\colon \mathbf{E}_p\longrightarrow \mathbf{E}'_p$ and
$(F_{\mathbf{A}})_p(v)=F_{\mathbf{A}}(v)=\mathbf{A}(p)(v)$ for every $v\in \mathbf{E}_p$. Since $\mathbf{A}_{F_{\mathbf{A}}}(p)$ and
$\mathbf{A}(p)$ agree on every vector in $\mathbf{E}_p$, we conclude that
$\mathbf{A}_{F_{\mathbf{A}}}(p)=\mathbf{A}(p)$ for every $p\in M$, and hence $\mathbf{A}_{F_{\mathbf{A}}}=\mathbf{A}$.

We conclude that the correspondences $F\mapsto \mathbf{A}_F$ and $\mathbf{A}\mapsto F_{\mathbf{A}}$ are inverse to each other, proving the bijection
\[
\underline{\operatorname{Hom}}(\mathbf{E},\mathbf{E}')\cong \Gamma\bigl(\operatorname{Hom}(\mathbf{E},\mathbf{E}')\bigr).
\]

(b) Define
\[
\Theta\colon \Gamma(\operatorname{Hom}(\mathbf{E},\mathbf{E}'))\longrightarrow \operatorname{Hom}_{C^{\infty}(M)}\bigl(\Gamma(\mathbf{E}),\Gamma(\mathbf{E}')\bigr),
\qquad
\Theta(\mathbf{A})(\boldsymbol{\sigma})(p):=\mathbf{A}(p)\bigl(\boldsymbol{\sigma}(p)\bigr).
\]
By Remark~\ref{obs F homomorfismo es C infty lineal}, $\Theta(\mathbf{A})$ is $C^{\infty}(M)$-linear.

We show that $\Theta$ is $C^{\infty}(M)$-linear in $\mathbf{A}$. Let $f\in C^{\infty}(M)$ and $\mathbf{A},\mathbf{B}\in \Gamma(\operatorname{Hom}(\mathbf{E},\mathbf{E}'))$. We have
\[
\Theta(f\mathbf{A}+\mathbf{B})(\boldsymbol{\sigma})(p)=(f\mathbf{A}+\mathbf{B})(p)(\boldsymbol{\sigma}(p))=(f(p)\mathbf{A}(p)+\mathbf{B}(p))(\boldsymbol{\sigma}(p))
\]
\[
=f(p)\mathbf{A}(p)(\boldsymbol{\sigma}(p))+\mathbf{B}(p)(\boldsymbol{\sigma}(p))=(f\Theta(\mathbf{A}))(\boldsymbol{\sigma})(p)+\Theta(\mathbf{B})(\boldsymbol{\sigma})(p),
\]
and therefore $\Theta(f\mathbf{A}+\mathbf{B})=f\Theta(\mathbf{A})+\Theta(\mathbf{B})$.

For surjectivity, let $\mathcal{F}\colon \Gamma(\mathbf{E})\longrightarrow\Gamma(\mathbf{E}')$ be a $C^{\infty}(M)$-linear map. By Lemma~\ref{lema de caracterizacion de homomorfismos de haces}, there exists a unique bundle homomorphism $F\colon \mathbf{E}\longrightarrow \mathbf{E}'$ such that $\mathcal{F}(\boldsymbol{\sigma})=F\circ\boldsymbol{\sigma}$ for every $\boldsymbol{\sigma}\in\Gamma(\mathbf{E})$. By (a), $F$ corresponds to a section $\mathbf{A}_F\in\Gamma(\operatorname{Hom}(\mathbf{E},\mathbf{E}'))$ and, by definition, $\Theta(\mathbf{A}_F)=\mathcal{F}$.

For injectivity, suppose that $\Theta(\mathbf{A}_1)=\Theta(\mathbf{A}_2)$. Then, for every section $\boldsymbol{\sigma}\in\Gamma(\mathbf{E})$ and every $p\in M$, we have $\mathbf{A}_1(p)\bigl(\boldsymbol{\sigma}(p)\bigr)=\mathbf{A}_2(p)\bigl(\boldsymbol{\sigma}(p)\bigr)$. Now let $v\in \mathbf{E}_p$. By Lemma~\ref{lema:seccion-valor-prescrito}, there exists a section $\boldsymbol{\sigma}\in\Gamma(\mathbf{E})$ with $\boldsymbol{\sigma}(p)=v$. Substituting into the preceding equality gives $\mathbf{A}_1(p)(v)=\mathbf{A}_2(p)(v)$. Since this holds for every $v\in \mathbf{E}_p$ and every $p\in M$, it follows that $\mathbf{A}_1=\mathbf{A}_2$.

Thus $\Theta$ is an isomorphism of $C^{\infty}(M)$-modules.

(c) The identification in part (c) is obtained by composing those in (a) and
(b). Compatibility with composition is justified as
follows. Let $F\colon \mathbf{E}\longrightarrow \mathbf{E}'$ and
$G\colon \mathbf{E}'\longrightarrow \mathbf{E}''$ be smooth bundle
homomorphisms over $M$. Their composition
$G\circ F\colon \mathbf{E}\longrightarrow \mathbf{E}''$ is also a
bundle homomorphism, and the associated section is
$\mathbf{A}_{G\circ F}(p)=(G\circ F)_p=G_p\circ F_p$ for $p\in M$.

On the other hand, the operator induced by $\mathbf{A}_{G\circ F}$ satisfies
\[
\Theta(\mathbf{A}_{G\circ F})(\boldsymbol{\sigma})(p)=\mathbf{A}_{G\circ F}(p)(\boldsymbol{\sigma}(p))=(G_p\circ F_p)(\boldsymbol{\sigma}(p)).
\]
Writing $\boldsymbol{\tau}=F\circ\boldsymbol{\sigma}
\in\Gamma(\mathbf{E}')$, we have
$\Theta(\mathbf{A}_{G\circ F})(\boldsymbol{\sigma})(p)
=G_p(\boldsymbol{\tau}(p))
=\Theta(\mathbf{A}_G)(\boldsymbol{\tau})(p)
=\Theta(\mathbf{A}_G)(F\circ\boldsymbol{\sigma})(p)$. Since
$F\circ\boldsymbol{\sigma}=\Theta(\mathbf{A}_F)(\boldsymbol{\sigma})$, it follows that
\[
\Theta(\mathbf{A}_{G\circ F})(\boldsymbol{\sigma})=\Theta(\mathbf{A}_G)(\Theta(\mathbf{A}_F)(\boldsymbol{\sigma}))=(\Theta(\mathbf{A}_G)\circ\Theta(\mathbf{A}_F))(\boldsymbol{\sigma}).
\]

Therefore, $\Theta(\mathbf{A}_{G\circ F})=\Theta(\mathbf{A}_G)\circ\Theta(\mathbf{A}_F)$, proving that the correspondence is compatible with composition of bundle homomorphisms.

Moreover, since (b) showed that $\Theta$ is an isomorphism of $C^\infty(M)$-modules, the identification is also compatible with the action of smooth functions; that is, $\Theta(\mathbf{A})(f\boldsymbol{\sigma})=f\cdot\Theta(\mathbf{A})(\boldsymbol{\sigma})$ for $f\in C^\infty(M)$ and $\boldsymbol{\sigma}\in\Gamma(\mathbf{E})$.

Consequently,
\[
\underline{\operatorname{Hom}}(\mathbf{E},\mathbf{E}')\cong \operatorname{Hom}_{C^\infty(M)}\bigl(\Gamma(\mathbf{E}),\Gamma(\mathbf{E}')\bigr)
\]
is a natural identification, compatible with both composition and the $C^\infty(M)$-module structure.
\end{proof}
The dual bundle is the special case of the homomorphism bundle in which the
target bundle is the trivial rank-one bundle.

\begin{definition}[Dual bundle]\label{def:haz-dual-suave}\index{dual bundle}
Let $M$ be a smooth manifold with or without boundary, and let $\mathbf E\to M$ be a
smooth real vector bundle. Its \emph{dual bundle} is
\[
 \mathbf E^*:=\operatorname{Hom}(\mathbf E,\underline{\mathbb R}_M).
\]
In particular, $(\mathbf E^*)_p=\operatorname{Hom}_{\mathbb R}
(\mathbf E_p,\mathbb R)$ for every $p\in M$.
\end{definition}

\begin{remark}[Dual frames and transition matrices]\label{remark:marcos-duales}
Let $\mathbf{E}\to M$ be a smooth vector bundle of rank $k$, and let $(\boldsymbol{\sigma}_1,\dots,\boldsymbol{\sigma}_k)$ be a smooth local frame of $\mathbf{E}$ defined on an open set $U\subseteq M$. The \emph{dual local frame} $(\boldsymbol{\sigma}^1,\dots,\boldsymbol{\sigma}^k)$ of $\mathbf{E}^*$ on $U$ is characterized by the condition
\[
\boldsymbol{\sigma}^i(\boldsymbol{\sigma}_j)=\delta^i_j,
\qquad 1\leq i,j\leq k.
\]

If $(\boldsymbol{\tau}_1,\dots,\boldsymbol{\tau}_k)$ is another smooth local frame of $\mathbf{E}$ on $U$, there exist smooth functions $A^j_i\colon U\longrightarrow\mathbb{R}$ such that
\[
\boldsymbol{\tau}_i=\displaystyle\sum_{j=1}^k A^j_i\boldsymbol{\sigma}_j,
\qquad i\in\{1,\dots,k\}.
\]
Thus the transition matrix between frames of $\mathbf{E}$ is $A=(A^j_i)\in GL(k,\mathbb{R})$.

The dual frame $(\boldsymbol{\tau}^1,\dots,\boldsymbol{\tau}^k)$ on $\mathbf{E}^*$ must satisfy
\[
\boldsymbol{\tau}^i(\boldsymbol{\tau}_j)=\delta^i_j,
\]
which implies, in index notation,
\[
\boldsymbol{\tau}^i=\displaystyle\sum_{j=1}^k (A^{-1})^i_j\boldsymbol{\sigma}^j.
\]

In matrix form, if $[\boldsymbol{\tau}]=A[\boldsymbol{\sigma}]$ denotes the relation between frames of $\mathbf{E}$, the dual frames satisfy
\[
[\boldsymbol{\tau}^T]=[\boldsymbol{\sigma}^T]A^{-1},
\]
so the transition matrix of $\mathbf{E}^*$ is
\[
B=(A^{-1})^T.
\]

Thus the transition matrices of the dual bundle $\mathbf{E}^*$ are the \emph{inverse transposes} of the transition matrices of $\mathbf{E}$.
\end{remark}

\paragraph{Conjugate spaces and bundles.}
In the complex case, we need to distinguish the ordinary dual from
conjugation. We first specify these constructions without imposing an
identification between a vector space and its dual.

\begin{definition}[Conjugate vector space]
\label{def:espacio-vectorial-conjugado}
Let $V$ be a complex vector space. Its \emph{conjugate space}
$\overline V$ is the same abelian group as $V$, with elements denoted by
$\overline v$ and scalar multiplication
\begin{equation}
\label{eq:producto-escalar-espacio-conjugado}
 \lambda\cdot\overline v
 :=\overline{\overline\lambda\,v}.
\end{equation}
The ordinary complex dual $V^*$ continues to mean
$\operatorname{Hom}_{\mathbb C}(V,\mathbb C)$.
\end{definition}

\begin{proposition}[Functorial properties of conjugation]
\label{prop:espacios-vectoriales-conjugados-funtoriales}
The preceding structure makes $\overline V$ a complex vector
space, and the map
\[
 c_V\colon V\longrightarrow\overline V,
 \qquad c_V(v)=\overline v,
\]
is additive and antilinear. The map
\[
 j_V\colon\overline{\overline V}\longrightarrow V,
 \qquad j_V(\overline{\overline v})=v,
\]
is a canonical complex-linear isomorphism.

If $A\colon V\to W$ is complex linear, there exists a
complex-linear map
\[
 \overline A\colon\overline V\longrightarrow\overline W,
 \qquad \overline A(\overline v)=\overline{Av}.
\]
We have
\[
 \overline{B\circ A}=\overline B\circ\overline A,
 \qquad
 \overline{\operatorname{id}_V}=\operatorname{id}_{\overline V},
 \qquad
 \overline{A^{-1}}=(\overline A)^{-1}
\]
when $A$ is invertible, and
$j_W\circ\overline{\overline A}=A\circ j_V$.
\end{proposition}

\begin{proof}
Addition is that of the underlying abelian group; that is,
$\overline v+\overline w=\overline{v+w}$. For
$\lambda,\mu\in\mathbb C$ and $v,w\in V$, the definition gives
\begin{align*}
 \lambda\cdot(\overline v+\overline w)
 &=\overline{\overline\lambda(v+w)}
  =\lambda\cdot\overline v+\lambda\cdot\overline w,\\
 (\lambda+\mu)\cdot\overline v
 &=\overline{(\overline\lambda+\overline\mu)v}
  =\lambda\cdot\overline v+\mu\cdot\overline v,\\
 \lambda\cdot(\mu\cdot\overline v)
 &=\overline{\overline\lambda\,\overline\mu\,v}
  =(\lambda\mu)\cdot\overline v,\\
 1\cdot\overline v&=\overline v.
\end{align*}
Together with the additive group axioms, these are all the vector
space identities. Moreover,
$c_V(v+w)=c_V(v)+c_V(w)$ and
$c_V(\lambda v)=\overline{\lambda v}
=\overline\lambda\cdot\overline v$, so $c_V$ is additive and
antilinear.
Under double conjugation,
\[
 j_V\bigl(\lambda\cdot\overline{\overline v}\bigr)
 =j_V\bigl(\overline{\,\overline\lambda
       \cdot_{\overline V}\overline v\,}\bigr)
 =j_V(\overline{\overline{\lambda v}})=\lambda v,
\]
which proves linearity of $j_V$; its inverse is
$v\mapsto\overline{\overline v}$.

For $\overline A$, if $\lambda\in\mathbb C$, then
\[
 \overline A(\lambda\overline v)
 =\overline{A(\overline\lambda v)}
 =\overline{\overline\lambda Av}
 =\lambda\overline{Av}.
\]
The remaining identities are checked on elements
$\overline v$, and hence hold as identities of maps.
\end{proof}

\begin{proposition}[Duality and conjugation]
\label{prop:dual-conjugado-identificacion-canonica}
There is a natural complex-linear isomorphism
\begin{equation}
\label{eq:kappa-dual-conjugado}
 \kappa_V\colon\overline{V^*}\longrightarrow(\overline V)^*,
 \qquad
 \kappa_V(\overline\ell)(\overline v)
 :=\overline{\ell(v)}.
\end{equation}
In particular, $V^*$, $\overline V$, $\overline{V^*}$, and
$(\overline V)^*$ are distinct objects; only
$\overline{V^*}\cong(\overline V)^*$ are canonically identified by
\eqref{eq:kappa-dual-conjugado}.
\end{proposition}

\begin{proof}
For $\mu\in\mathbb C$,
\[
 \kappa_V(\overline\ell)(\mu\overline v)
 =\overline{\ell(\overline\mu v)}
 =\mu\,\overline{\ell(v)},
\]
so $\kappa_V(\overline\ell)\in(\overline V)^*$. Likewise,
\[
 \kappa_V(\lambda\overline\ell)(\overline v)
 =\kappa_V(\overline{\overline\lambda\ell})(\overline v)
 =\lambda\,\overline{\ell(v)},
\]
so $\kappa_V$ is complex linear. The formula
\[
 \ell(v):=\overline{L(\overline v)},
 \qquad L\in(\overline V)^*,
\]
defines a complex-linear functional, since
$\overline{\lambda v}=\overline\lambda\cdot\overline v$ in
$\overline V$, and is the two-sided inverse of \eqref{eq:kappa-dual-conjugado}.
Finally, if $A\colon V\to W$ is complex linear and
$A^*\ell=\ell\circ A$, then
\[
 \kappa_V\circ\overline{A^*}
 =(\overline A)^*\circ\kappa_W,
\]
as one checks by evaluating both sides at $\overline v$; this proves the
asserted naturality.
\end{proof}

For complex bundles, homomorphisms are complex linear on each
fiber, in accordance with
Definition~\ref{def:nociones-fundamentales-homomorfismo-de-haces}. The construction of
$\operatorname{Hom}(\mathbf E,\mathbf F)$ in
Proposition~\ref{haz Hom vectorial} applies unchanged, replacing
$\mathbb R$ by $\mathbb C$ and requiring complex linearity. Likewise,
for a complex bundle we set
\[
 \mathbf E^*:=\operatorname{Hom}_{\mathbb C}
 (\mathbf E,\underline{\mathbb C}_M).
\]

\begin{definition}[Conjugate bundle]
\label{def:haz-vectorial-conjugado}
Let $M$ be a smooth manifold with or without boundary and let $\pi\colon \mathbf{E}\to M$ be a smooth complex vector bundle. The
\emph{conjugate bundle} $\overline{\mathbf{E}}\to M$ has the same total space,
underlying smooth real manifold, projection, and fiberwise addition
as $\mathbf{E}$. We write $\overline v$ when an element $v\in \mathbf{E}$ is
regarded as an element of $\overline{\mathbf{E}}$, and define scalar multiplication by
\begin{equation}
\label{eq:producto-escalar-haz-conjugado-intrinseco}
 \lambda\cdot_{\overline{\mathbf{E}}}\overline v
 :=\overline{\overline\lambda\,v},
 \qquad \lambda\in\mathbb C.
\end{equation}
\end{definition}

\begin{proposition}[Construction and naturality of the conjugate bundle]
\label{prop:conjugacion-haces-secciones-morfismos}
Let $M$ be a smooth manifold with or without boundary and let
$\mathbf{E}\longrightarrow M$ be a smooth complex vector bundle. Each fiber
$\overline{\mathbf{E}}_p$ is the conjugate vector space of $\mathbf{E}_p$, and the
projection $\overline{\mathbf{E}}\to M$ is smooth. If
\[
 \Phi_\alpha\colon \mathbf{E}|_{U_\alpha}\longrightarrow
 U_\alpha\times\mathbb C^r,
 \qquad \Phi_\alpha(v)=(p,z),
\]
is a smooth complex trivialization of $\mathbf{E}$, then
\[
 \overline\Phi_\alpha\colon\overline{\mathbf{E}}|_{U_\alpha}
 \longrightarrow U_\alpha\times\mathbb C^r,
 \qquad
 \overline\Phi_\alpha(\overline v):=(p,\overline z),
\]
is a smooth complex trivialization of $\overline{\mathbf{E}}$. If the transition
functions of $\mathbf{E}$ are given by
\[
 (\Phi_\beta\circ\Phi_\alpha^{-1})(p,z)
 =(p,g_{\beta\alpha}(p)z),
\]
the transition functions computed for $\overline{\mathbf{E}}$ are
$\overline{g_{\beta\alpha}}$. Consequently,
\[
 \overline{g_{\gamma\beta}}\,
 \overline{g_{\beta\alpha}}
 =\overline{g_{\gamma\alpha}},
 \qquad
 \overline{g_{\alpha\beta}}
 =(\overline{g_{\beta\alpha}})^{-1}.
\]
The construction is independent of the chosen trivializations, since
it is defined on the same total space by the intrinsic operation
\eqref{eq:producto-escalar-haz-conjugado-intrinseco} on each fiber.
The assignment on sections
\[
 \Gamma(\mathbf{E})\longrightarrow\Gamma(\overline{\mathbf{E}}),
 \qquad \mathbf{u}\longmapsto\overline{\mathbf{u}},
 \qquad \overline{\mathbf{u}}(x):=\overline{\mathbf{u}(x)},
\]
is an additive, antilinear bijection compatible with restrictions.

If $A\colon \mathbf{E}\to \mathbf{F}$ is a homomorphism of complex bundles, then
\[
 \overline A\colon\overline{\mathbf{E}}\longrightarrow\overline{\mathbf{F}},
 \qquad \overline A(\overline v)=\overline{A(v)},
\]
is a smooth homomorphism of complex bundles. Conjugation preserves identities,
compositions, and inverses, and there is a natural identification
$j_{\mathbf{E}}\colon\overline{\overline{\mathbf{E}}}\xrightarrow{\cong}\mathbf{E}$. Finally,
\[
 \kappa_{\mathbf{E}}\colon\overline{\mathbf{E}^*}\longrightarrow(\overline{\mathbf{E}})^*,
 \qquad
 \kappa_{\mathbf{E}}(\overline\ell)(\overline v)=\overline{\ell(v)},
\]
is a smooth complex-linear bundle isomorphism.
\end{proposition}

\begin{proof}
Proposition~\ref{prop:espacios-vectoriales-conjugados-funtoriales},
applied fiberwise, proves the complex vector space axioms on
$\overline{\mathbf{E}}_p$. The underlying real manifold and the projection are
unchanged; therefore, $\overline{\mathbf{E}}\to M$ remains a smooth map.

Let $\Phi_\alpha(v)=(p,z)$. Since conjugation
$z\mapsto\overline z$ is a real-linear isomorphism of $\mathbb C^r$, the
map $\overline\Phi_\alpha$ is a real diffeomorphism. If
$\lambda\in\mathbb C$, then
\[
 \overline\Phi_\alpha
 (\lambda\cdot_{\overline{\mathbf{E}}}\overline v)
 =\overline\Phi_\alpha(\overline{\overline\lambda v})
 =(p,\lambda\overline z),
\]
so it is complex linear on each fiber. Thus these functions are
smooth complex trivializations of the structure already defined.

To compute their transition maps, start with a vector whose coordinate is $w$
in the conjugate trivialization $\overline\Phi_\alpha$. The same element,
regarded as an element of $\mathbf{E}$, has coordinate $\overline w$; hence
\[
 (\overline\Phi_\beta\circ\overline\Phi_\alpha^{-1})(p,w)
 =\bigl(p,\overline{g_{\beta\alpha}(p)\overline w}\bigr)
 =\bigl(p,\overline{g_{\beta\alpha}(p)}w\bigr).
\]
The two identities in the statement now follow by conjugating the
elementary identities already satisfied by the transition functions of
$\mathbf{E}$. These were not used to construct $\overline{\mathbf{E}}$: they are consequences of
the intrinsic construction. Since the total space, the projection, and
scalar multiplication
\eqref{eq:producto-escalar-haz-conjugado-intrinseco} are all independent of any chart,
independence of the trivializations is immediate.

If $u_\alpha$ are the local components of $\mathbf{u}$, then
$u_\beta=g_{\beta\alpha}u_\alpha$. Therefore,
$\overline{u_\beta}=\overline{g_{\beta\alpha}}\,
\overline{u_\alpha}$, which are exactly the transformation rules for a
section of $\overline{\mathbf{E}}$. The reverse argument uses conjugation once more.
The scalar multiplication rule on $\overline{\mathbf{E}}$ proves
$\overline{f\mathbf{u}}=\overline f\,\overline{\mathbf{u}}$ for
$f\in C^\infty(M,\mathbb C)$; compatibility with restrictions is
pointwise.

In local frames, if $A$ has matrix $A_\alpha(x)$, then
$\overline A$ has matrix $\overline{A_\alpha(x)}$; the change-of-frame rules
show that these matrices glue, and smoothness is immediate.
Complex linearity of $\overline A$ is checked using
\eqref{eq:producto-escalar-haz-conjugado-intrinseco}. The functorial
properties and the formula
$j_{\mathbf{E}}(\overline{\overline v})=v$ are verified fiberwise by
Proposition~\ref{prop:espacios-vectoriales-conjugados-funtoriales}; in
conjugate trivializations, $j_{\mathbf{E}}$ has the identity matrix and is smooth. The last
assertion is the bundle version of
\eqref{eq:kappa-dual-conjugado}; its local formula and that of its inverse have
constant matrices, so both maps are smooth.
\end{proof}

\begin{corollary}[Transition functions of dual and homomorphism bundles]
\label{cor:cociclos-hom-dual-haces-complejos}
Let $M$ be a smooth manifold with or without boundary, and let
$\mathbf{E},\mathbf{F}\longrightarrow M$ be smooth complex vector bundles.
With the preceding convention, the transition functions of $\mathbf{E}^*$ are
$(g_{\beta\alpha}^{-1})^T$, those of $\overline{\mathbf{E}}$ are
$\overline{g_{\beta\alpha}}$, and those of
$\operatorname{Hom}(\mathbf{E},\mathbf{F})$ act on a matrix $A$ by
\[
 A\longmapsto h_{\beta\alpha}A g_{\beta\alpha}^{-1}.
\]
Conjugating them gives precisely the transition functions of the corresponding
conjugate bundles. These formulas do not identify $\mathbf{E}$ with $\mathbf{E}^*$.
\end{corollary}

\begin{proof}
The first formula follows by requiring the evaluation
$\mathbf{E}^*\otimes \mathbf{E}\to\underline{\mathbb C}$ to be independent of the frame. The
second is the definition, and the third follows from
$A_\beta(g_{\beta\alpha}v)=h_{\beta\alpha}A_\alpha(v)$. Conjugating each
equality gives the last assertion.
\end{proof}

 Certain sections of the tangent bundle $TM$ receive a special name because of their importance in the study of smooth manifolds. These sections are called vector fields:
 \begin{definition}\label{def:nociones-fundamentales-campo-vectorial-en}\index{vector field}
 Let $M$ be a smooth manifold with or without boundary. A \textbf{vector
 field on $M$} is a section
 $\mathbf{X}\colon M\longrightarrow TM$ of the map
 $\pi\colon TM\longrightarrow M$. We usually write
 $\mathbf{X}(p)$ as $\mathbf{X}_{p}$ for each $p\in M$.
 \end{definition}
\begin{note}
 The set of smooth vector fields defined on a smooth manifold $M$ with or without boundary is usually denoted by $\mathfrak{X}(M)$.
\end{note}
 By Proposition~\ref{basedetpm}, every vector field
 $\mathbf{X}\colon M\longrightarrow TM$ can be written in terms of
 a smooth chart $(U,\phi)$ with $\phi=(x^{1},\dots,x^{n})$ as
 $\mathbf{X}=X^{i}\boldsymbol{\partial}_{i}$, where,
 for each $i\in \{1,\dots,n\}$,
 $\boldsymbol{\partial}_{i}\colon U\longrightarrow TU$
 is the vector field defined on $U$ by
 $\boldsymbol{\partial}_{i}(p)
 :=\boldsymbol{\partial}_{i}\biggr|_{p}$ for each
 $p\in U$. The functions $X^{i}$ are called the
 \textit{\textbf{component functions of $\mathbf{X}$ on $U$}}.

 We have the following criterion for smoothness of a vector field:

 \begin{proposition}[Regularity of a vector field in coordinates]
\label{criterio para la suavidad de un campo vectorial}
Let $M$ be a smooth manifold with or without boundary and let
$\mathbf X\colon M\to TM$ be a raw vector field.
For $k\in\mathbb N_0\cup\{\infty\}$ and any smooth chart
$(U,(x^i))$, the field $\mathbf X\restriction_U$ is of class $C^k$
if and only if its component functions $X^i$ are of class $C^k$.
In particular, the field is continuous if and only if its components
are continuous in the charts of an atlas covering $M$.
\end{proposition}
\begin{proof}
The coordinate fields $(\boldsymbol\partial_i)_{i=1}^n$ form a smooth
frame of the tangent bundle over $U$. In its associated trivialization,
the field is represented by $p\mapsto(p,(X^1(p),\dots,X^n(p)))$.
Proposition~\ref{suavidad de secciones en bruto locales} gives the
criterion for each $k$, including continuity when $k=0$.
\end{proof}

The vector field of vector calculus is generalized by assigning to each
$p\in M$ a vector $\mathbf X_p\in T_pM$. Vectors at different points
generally belong to different vector spaces; the tangent bundle brings
them together so that continuity of the assignment can be considered.
When $M$ is embedded in $\mathbb R^N$, they can be drawn as ambient
arrows based at $p$, but each arrow must belong to the corresponding
tangent plane.

In charts $x=\phi(p)$ and $y=\psi(p)$, setting
$T:=\psi\circ\phi^{-1}$ and
$\widehat X_\phi:= (X^1\circ\phi^{-1},\dots,X^n\circ\phi^{-1})$, the
chain rule gives
\begin{equation}\label{eq:compatibilidad-componentes-campos}
 \widehat X_\psi(T(x))=DT(x)\widehat X_\phi(x).
\end{equation}
Indeed, component $j$ in the $y$ coordinates is
$\mathbf X_p(y^j)=\displaystyle\sum_{i=1}^n
 X^i(p)\displaystyle\frac{\partial(y^j\circ\phi^{-1})}{\partial x^i}(\phi(p))$.
This is the condition for defining a field from formulas in different
charts. In addition to the change of argument occurring for scalar
functions, the vector is transformed here by the Jacobian.

\begin{example}[A smooth rotation field on the sphere]
\label{ej:campo-rotacion-esfera-autocontenido}
Use the atlas of Example~\ref{ej:altura-esfera-cartas-autocontenido}
and the inclusion of $S^2$ in $\mathbb R^3$. Under the tangent
identification induced by that inclusion,
\[
 T_pS^2=\{w\in\mathbb R^3\mid\langle p,w\rangle=0\}.
\]
To verify the formula directly, differentiating $|\Phi_N(u,v)|^2=1$
gives $\langle\Phi_N,D\Phi_N(\xi)\rangle=0$ for all $\xi\in\mathbb R^2$.
The differential of $\Phi_N$ is injective: the identity
$\phi_N\circ\Phi_N=\operatorname{Id}$, with $\phi_N$ extended by its
formula to the ambient open set $z\neq1$, gives a left inverse for its
differential. Thus its image has dimension two and equals the kernel of
the nonzero functional $w\mapsto\langle p,w\rangle$.
The same proof using $\Phi_S$ covers the remaining pole.

Define
\[
 \mathbf X_p:=(-y,x,0),\qquad p=(x,y,z)\in S^2.
\]
One has $\langle p,\mathbf X_p\rangle=-xy+yx=0$, so the arrows are
tangent. By contrast, the ambient vector $p$ is normal and does not
belong to $T_pS^2$, since $\langle p,p\rangle=1$.

To compute the components of $\mathbf X$, use
\[
 c_p(t):=(x\cos t-y\sin t,\ x\sin t+y\cos t,\ z).
\]
This curve remains in $S^2$, equals $p$ at $t=0$, and has velocity
$\mathbf X_p$ there. In both charts,
\[
 \phi_\alpha(c_p(t))
 =(u\cos t-v\sin t,\ u\sin t+v\cos t),
 \qquad \alpha\in\{N,S\}.
\]
Thus the component column is $(-v,u)$, and
\[
 \mathbf X=-v\,\boldsymbol\partial_u+u\,\boldsymbol\partial_v
\]
in each of the two charts. The coefficients are polynomial, so the
field is smooth, including at the poles, where both arrows vanish.

Compatibility can be checked without using the ambient definition.
For $q=(u,v)\neq0$, the stereographic transition is $T(q)=r^{-2}q$, and
\[
 DT(q)=r^{-2}I-2r^{-4}q q^{\mathsf T}.
\]
Since $\langle q,(-v,u)\rangle=0$,
\[
 DT(q)(-v,u)=r^{-2}(-v,u)
 =\widehat X_S(T(q)).
\]
Thus~\eqref{eq:compatibilidad-componentes-campos} holds.
The ambient length of the arrow is $\sqrt{x^2+y^2}=\sqrt{1-z^2}$:
it is maximal at the equator and zero at the poles.
The curves $c_p$ preserve height and run along horizontal circles
in the direction of the field.
\end{example}

\begin{figura}[htbp]
 \includegraphics[width=0.88\linewidth]{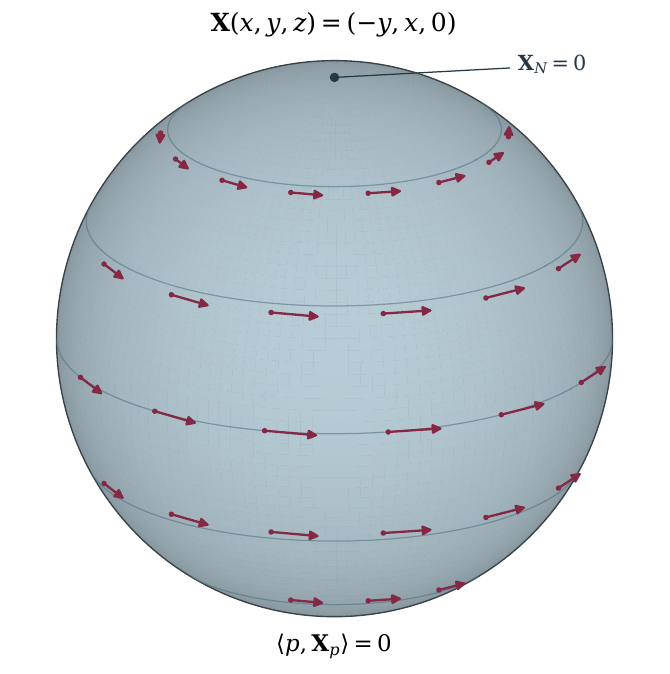}
 \caption{The tangent field $\mathbf X_p=(-y,x,0)$ on $S^2$.
 Arrows start at the indicated points and represent vectors in their
 tangent planes; the same length scale is used for every arrow.
 The circles show rotation trajectories. The field vanishes at the poles,
 although both the sphere and the field remain smooth.}
 \label{fig:campo-vectorial-esfera-rotacion}
\end{figura}

\begin{example}[Fields of lower regularity]
\label{ej:campos-continuos-y-discontinuos-esfera}
With $h$ and $\mathbf X$ as above, the field $\mathbf Y:=|h|\mathbf X$
is continuous. It is not $C^1$ at $p=(1,0,0)$: along
$c(t)=(\cos t,0,\sin t)$, its second ambient component is
$|\sin t|\cos t$, whose difference quotients at $0$ have limits $1$
and $-1$. If $\mathbf Y$ were $C^1$, this component of its ambient
representation along a smooth curve would also be $C^1$.
The raw field $\mathbf Z:=\mathbf 1_{\{h>0\}}\mathbf X$ is discontinuous
at the equator, since it tends to $\mathbf X_p\neq0$ from one side and
to $0$ from the other. Both examples are tangent at every point;
tangency and regularity are different conditions.
\end{example}
 The preceding proposition tells us that the coordinate vector fields $\boldsymbol{\partial}_{i}$ are smooth vector fields on $U$.

 An essential property of vector fields is that they define operators on spaces of smooth real-valued functions. If $\mathbf{X}\in\mathfrak{X}(M)$, $U\subseteq M$ is open, and $f\in C^{\infty}(U)$, we obtain a new function $\mathbf{X}(f)\colon U\longrightarrow \mathbb{R}$ defined by \[(\mathbf{X}(f))(p)=\mathbf{X}_{p}(f).\]

This terminology gives another smoothness criterion for vector fields.

\begin{proposition}\label{prop: criterios de suavidad para un campo vectorial suave}
 Let $M$ be a smooth manifold with or without boundary, and let
 $\mathbf{X}\colon M\longrightarrow TM$ be a rough vector field. The
 following are equivalent:
 \begin{enumerate}[label=(\alph*)]
 \item $\mathbf{X}$ is smooth.
 \item For each $f\in C^{\infty}(M)$, the function $\mathbf{X}(f)$ is smooth on $M$.
 \item For every open subset $U\subseteq M$ and every $f\in C^{\infty}(U)$, the function $(\mathbf{X}\restriction_U)(f)$ is smooth on $U$.
 \end{enumerate}
\end{proposition}
\begin{proof}
If \(\mathbf X=X^i\boldsymbol\partial_i\) is smooth in a chart \(U\), then for
\(f\in C^\infty(U)\) we have
\[
(\mathbf X\restriction_U)(f)
=\displaystyle\sum_{i=1}^{n}X^i\frac{\partial f}{\partial x^i}\in C^\infty(U).
\]
This proves \((a)\Rightarrow(c)\), and \((c)\Rightarrow(b)\) follows
by taking \(U=M\).

Suppose (b) holds. Fix \(p\in M\), a chart \((U,x)\) around \(p\),
and an open set \(W\) with \(p\in W\Subset U\). Let
\(\chi\in C^\infty(M)\) equal \(1\) on \(\overline W\) and have support
contained in \(U\). The function defined by \(\chi x^i\) on \(U\) and by
zero outside \(U\) is smooth on \(M\). On \(W\),
\[
\mathbf X(\chi x^i)=\mathbf X(x^i)=X^i.
\]
Assumption (b) shows that each component \(X^i\) is smooth near
\(p\). Since \(p\) was arbitrary,
Proposition~\ref{criterio para la suavidad de un campo vectorial} proves
\((b)\Rightarrow(a)\).
\end{proof}
 Finally, we introduce the notions of a cotangent vector, cotangent space, and cotangent bundle of a smooth manifold $M$ with or without boundary. We have the following definitions:
 \begin{definition}\label{def:nociones-fundamentales-espacio-cotangente-a-en}\index{cotangent space}
 Let $M$ be a smooth manifold with or without boundary. For each $p\in M$, define the \textbf{cotangent space to $M$ at $p$}, denoted by $T_{p}^{*}M$, as the dual space of $T_{p}M$, \[T_{p}^{*}M:=(T_{p}M)^{*}.\]
 \end{definition}
 \begin{remark}\label{obs:nociones-fundamentales-campo-covectorial-coordenado}
 For each $p\in M$, if $(U,\phi)$ is a coordinate chart on an open subset $U\subseteq M$ with $p\in U$, the basis $\left(\boldsymbol{\partial}_{i}\biggr|_{p}\right)$ gives a dual basis for $T_{p}^{*}M$. Denote this basis by $(\boldsymbol{\lambda}^{i}|_{p})$. Now recall the differential of a function and that, if $f\in C^{\infty}(M)$, then $T_{f(p)}\mathbb{R}\cong \mathbb{R}$. We may therefore regard $df_{p}(v)\in \mathbb{R}$ for each $v\in T_{p}M$ and each $p\in M$. Then $df_{p}\in T_{p}^{*}M$ for each $f\in C^{\infty}(M)$. We know that $df_{p}=A_{i}(p)\boldsymbol{\lambda}^{i}|_{p}$ for some function $A_{i}\colon U\longrightarrow \mathbb{R}$.

 By the definition of the differential and the fact that $(\boldsymbol{\lambda}^{i}|_{p})$ is the dual basis of $\boldsymbol{\partial}_{i}\biggr|_{p}$, \[A_{i}(p)=df_{p}\left(\boldsymbol{\partial}_{i}\biggr|_{p}\right)=\boldsymbol{\partial}_{i}\biggr|_{p}(f)=\frac{\partial f}{\partial x^{i}}(p)\]

 so
 \[df_{p}=\frac{\partial f}{\partial x^{i}}(p)\boldsymbol{\lambda}^{i}|_{p}\]

 In particular, taking $x^{j}\in C^{\infty}(U)$, for each $p\in U$ we have $\mathbf{d}x^{j}_{p}\in T_{p}^{*}(M)$, so \[\mathbf{d}x^{j}|_{p}=\frac{\partial x^{j}}{\partial x^{i}}(p)\boldsymbol{\lambda}^{i}|_{p}=\delta^{j}_{i}\boldsymbol{\lambda}^{i}|_{p}=\boldsymbol{\lambda}^{j}|_{p}\] Thus the dual basis of $T_{p}^{*}(M)$ is in fact $(\mathbf{d}x^{i}|_{p})$. The assignment $p\mapsto \mathbf{d}x^{i}|_{p}$ defines a covector field (for now, a rough one) on $U$ called a \textit{\textbf{coordinate covector field}}.
 \end{remark}
 We define the cotangent bundle:
 \begin{definition}\label{def:nociones-fundamentales-haz-cotangente-de}\index{cotangent bundle}
 If $M$ is a smooth manifold with or without boundary, define the \textbf{cotangent bundle of $M$} by \[T^{*}M:=\coprod_{p\in M}T_{p}^{*}M.\]
 \end{definition}
 Of course, for the cotangent bundle we can define an analogue of vector fields on the tangent bundle.
 \begin{definition}\label{def:nociones-fundamentales-campo-covectorial}\index{covector field}
 Let $M$ be a smooth manifold with or without boundary. A local or global section of $T^{*}M$ is called a
\textbf{covector field}. A \textbf{rough covector field} is a
local or global section that need not be continuous, and a
\textbf{smooth covector field} is a smooth local or global section.
 \end{definition}
 The smooth structure is obtained by dualizing the changes of basis in the tangent
 bundle. \begin{proposition}\label{prop:nociones-fundamentales-variedad-suave-frontera-dimension-proyeccion-estandar}
 Let $M$ be a smooth manifold with or without boundary of dimension $n$. With the standard projection $\pi\colon T^{*}M\longrightarrow M$ such that $\pi(\omega_{p})=p$ if $\omega_{p}\in T_{p}^{*}M$, and with the vector space structure on each fiber of $\pi$ inherited from that of $T_{p}M$, the cotangent bundle $T^{*}M$ has a unique topology and a unique smooth structure making it a smooth vector bundle of rank $n$ over $M$ for which the coordinate covector fields are smooth local sections. Moreover, if $U\subseteq M$ is the domain of a smooth chart with coordinates $\phi=(x^{1},\dots,x^{n})$, the map \[\widetilde{\phi}(\xi_{i}\mathbf{d}x^{i}|_{p})=(x^{1}(p),\dots,x^{n}(p),\xi_{1},\dots,\xi_{n})\] is a smooth coordinate chart for $T^{*}M$.
 \end{proposition}
\begin{proof}
For a chart $(U,x)$, with $x=(x^1,\ldots,x^n)$, define the
trivialization
\[
 \Phi_x\colon\pi^{-1}(U)\longrightarrow U\times\mathbb R^n,
 \qquad
 \Phi_x\left(\displaystyle\sum_{i=1}^n\xi_i\,\mathbf{d}x^i|_p\right)
 =(p,(\xi_1,\ldots,\xi_n)).
\]
It is a fiber-preserving bijection and is linear on each fiber. Let $(U,x)$ and
$(V,y)$ be two charts, and set
\[
 J(p):=D(y\circ x^{-1})(x(p)).
\]
The change-of-basis formula
\eqref{eq:cambio-base-vectores-coordenados} says that the components of a
vector satisfy $v_y=Jv_x$. The identity
$\omega(v)=\xi_x^Tv_x=\xi_y^Tv_y$ for every $v$ then implies
\[
 \xi_y=(J^{-1})^T\xi_x.
\]
Consequently, the change of trivialization of the cotangent bundle is
\[
 (\Phi_y\circ\Phi_x^{-1})(p,\xi)
 =\bigl(p,(J(p)^{-1})^T\xi\bigr).
\]
The map $p\mapsto(J(p)^{-1})^T$ is smooth and takes values in
$GL(n,\mathbb R)$. Lemma~\ref{lema del haz vectorial suave} gives
a unique topology and a unique smooth structure for which the maps
$\Phi_x$ are local trivializations. In a chart of the total space, the
projection is $(x,\xi)\mapsto x$ and the fields $\mathbf{d}x^i$ have constant
components; hence the projection and the coordinate covector
fields are smooth. The same calculation works on relatively open subsets
of $\mathbb H^n$ and covers the boundary case.
\end{proof}

 We have the following smoothness criterion for a covector field:
 \begin{proposition}\label{prop:nociones-fundamentales-variedad-suave-frontera-campo-covectorial-bruto}
 Let $M$ be a smooth manifold with or without boundary and let
 $\boldsymbol{\omega}\colon M\longrightarrow T^{*}M$ be a rough
 covector field. The following statements are equivalent:
 \begin{enumerate}[label=(\alph*)]
 \item $\boldsymbol{\omega}$ is smooth.
 \item In every coordinate chart, the component functions of $\boldsymbol{\omega}$ are smooth.
 \item Every point of $M$ lies in some coordinate chart in which $\boldsymbol{\omega}$ has smooth component functions.
 \item For every smooth vector field $\mathbf{X}\in \mathfrak{X}(M)$, the function $\boldsymbol{\omega}(\mathbf{X})$ is smooth on $M$.
 \item For every open subset $U\subseteq M$ and every smooth vector field $\mathbf{X}$ on $U$, the function $\boldsymbol{\omega}(\mathbf{X})\colon U\longrightarrow \mathbb{R}$ is smooth on $U$.
 \end{enumerate}
 \end{proposition}
\begin{proof}
In a chart \((U,x)\), write
\(\boldsymbol\omega=\omega_i\,\mathbf d x^i\). The component criterion
in Proposition~\ref{suavidad de secciones en bruto locales} proves the
equivalence of (a), (b), and (c).

If (a) holds and
\(\mathbf X=X^i\boldsymbol\partial_i\) is smooth on \(U\), then
\(\boldsymbol\omega(\mathbf X)=\omega_iX^i\) is smooth. Therefore,
\((a)\Rightarrow(e)\Rightarrow(d)\).

Suppose (d) holds. For a chart \((U,x)\), a point \(p\in U\), and an open set
\(W\) with \(p\in W\Subset U\), choose
\(\chi\in C^\infty(M)\) equal to \(1\) on \(\overline W\) and supported in
\(U\). The field defined by
\(\chi\boldsymbol\partial_i\) on \(U\) and by zero outside \(U\) is smooth
on \(M\). On \(W\),
\[
\boldsymbol\omega(\chi\boldsymbol\partial_i)=\omega_i.
\]
Assumption (d) shows that the components are smooth near every
point, proving (c).
\end{proof}

 This has the following consequence:
 \begin{corollary}\label{cor:nociones-fundamentales-diferencial-funcion-suave-campo-covectorial-suave}
 Let $M$ be a smooth manifold with or without boundary. The differential of every
 function $f\in C^\infty(M)$ is a smooth covector field.
 \end{corollary}
\begin{proof}
In coordinates,
\[
df=\displaystyle\sum_{i=1}^{n}\frac{\partial f}{\partial x^i}\,\mathbf d x^i.
\]
Its components are smooth, so the preceding proposition applies.
\end{proof}
 Finally, before concluding this section, we introduce the concept of a smooth vector subbundle:
 \begin{definition}\label{def:nociones-fundamentales-subhaz-vectorial-de}\index{vector subbundle}
 If $\pi_{\mathbf{E}}\colon \mathbf{E}\longrightarrow M$ is a vector
 bundle, a \textbf{vector subbundle of $\mathbf{E}$} is a vector
 bundle $\pi_{\mathbf{D}}\colon \mathbf{D}\longrightarrow M$, where
 $\mathbf{D}\subseteq \mathbf{E}$ carries the subspace topology and
 $\pi_{\mathbf{D}}$ is the restriction of $\pi_{\mathbf{E}}$ to
 $\mathbf{D}$, such that, for each $p\in M$,
 $\mathbf{D}_{p}=\mathbf{D}\cap \mathbf{E}_{p}$ is a vector
 subspace of $\mathbf{E}_{p}$ and all the fibers $\mathbf{D}_{p}$ have
 the same dimension. If $M$ is a smooth manifold with or without boundary and
 $\mathbf{E}\longrightarrow M$ is a smooth vector bundle, a
 \textbf{smooth vector subbundle of $\mathbf{E}$} is a smooth vector
 bundle that is a vector subbundle of $\mathbf{E}$ and an embedded
 submanifold of $\mathbf{E}$.
 \end{definition}
\begin{figura}[h]
    \includegraphics[scale=1]{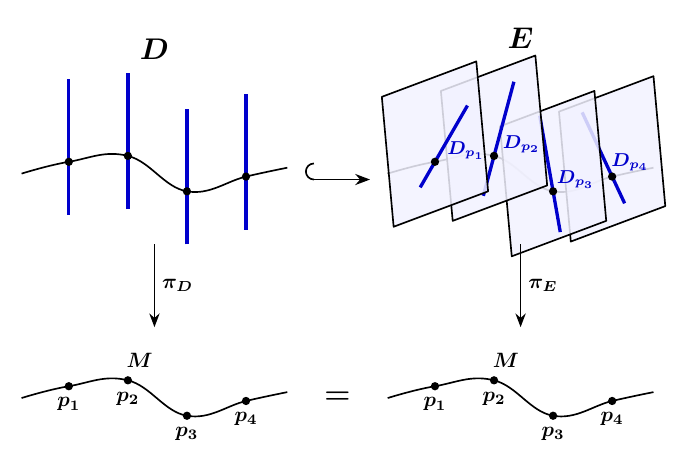}
    \caption{Smooth vector subbundle.}
    \label{fig:subhaz-vectorial-suave}
\end{figura}
 \begin{lemma}[Smooth local frame criterion for subbundles]\label{criterio del marco local suave para subhaces}\index{smooth local frame criterion for subbundles}
 Let $M$ be a smooth manifold with or without boundary, and let
 $\pi\colon \mathbf{E}\longrightarrow M$ be a smooth vector bundle.
 Suppose that, for each $p\in M$, we are given a vector subspace of
 fixed finite dimension $m$, $\mathbf{D}_{p}\subseteq \mathbf{E}_{p}$.
 Then $\mathbf{D}=\displaystyle\bigcup_{p\in M}\mathbf{D}_{p}
 \subseteq \mathbf{E}$ is a smooth vector subbundle of $\mathbf{E}$ if and
 only if every point of $M$ has a neighborhood $U$ on which there exist
 smooth local sections
 $\boldsymbol{\sigma}_{1},\dots,\boldsymbol{\sigma}_{m}\colon
 U\longrightarrow \mathbf{E}$ such that
 $\boldsymbol{\sigma}_{1}(q),\dots,\boldsymbol{\sigma}_{m}(q)$ form a
 basis for $\mathbf{D}_{q}$ for each $q\in U$.

 \end{lemma}
\begin{proof}
If \(\mathbf D\) is a smooth vector subbundle, a local trivialization of
\(\mathbf D\) gives a smooth local frame
\(\boldsymbol\sigma_1,\dots,\boldsymbol\sigma_m\). The smooth inclusion
\(\mathbf D\hookrightarrow\mathbf E\) allows us to regard its elements
as smooth sections of \(\mathbf E\), and their values form a basis for
\(\mathbf D_q\).

Conversely, suppose such sections exist around
\(p\in M\). Choose a local frame
\(\boldsymbol e_1,\dots,\boldsymbol e_r\) of \(\mathbf E\). In the fiber
over \(p\), extend
\(\boldsymbol\sigma_1(p),\dots,\boldsymbol\sigma_m(p)\) to a basis for
\(\mathbf E_p\) using some of the vectors
\(\boldsymbol e_j(p)\). The determinant of the matrix of these vectors in
the frame \((\boldsymbol e_j)\) is nonzero at \(p\); by continuity,
it remains nonzero after shrinking the neighborhood. We thus obtain
a local frame of \(\mathbf E\)
\[
\boldsymbol\sigma_1,\dots,\boldsymbol\sigma_m,
\boldsymbol\tau_{m+1},\dots,\boldsymbol\tau_r.
\]
In the trivialization determined by this frame,
\[
\mathbf D\restriction_U
\longleftrightarrow U\times(\mathbb R^m\times\{0\})
\subseteq U\times\mathbb R^r.
\]
Therefore, \(\mathbf D\restriction_U\to U\) is a smooth vector bundle of
rank \(m\), and its total space is an embedded submanifold of
\(\mathbf E\restriction_U\). On intersections, changes of frame
between the first \(m\) sections take values in \(GL(m,\mathbb R)\) and
are smooth, so these local structures glue. Thus
\(\mathbf D\) is a smooth vector subbundle of \(\mathbf E\).
\end{proof}
 \section{Integral curves and flows}

A vector field can be read as a law of motion: at each point, it prescribes a velocity. The trajectories obeying this law are its integral curves. Local existence and uniqueness for ordinary differential equations allow us to assemble them into a family of diffeomorphisms between open subsets of the manifold, called a flow.

The existence time may depend on the initial point, so a flow need not be defined for every $t\in \mathbb{R}$ or on the entire manifold. When it is, it determines a smooth action of $\mathbb{R}$ on $M$.

The next definition specifies this relation between velocity and vector field.

\begin{definition}\label{def:nociones-fundamentales-integral-variedad-suave-frontera-curva}\index{integral curve}\index{initial point!of an integral curve}
Let $M$ be a smooth manifold without boundary and let $\mathbf{X}\in \mathfrak{X}(M)$. We say that a smooth curve $\gamma\colon I\longrightarrow M$ is an integral curve of $\mathbf{X}$ if, for every $t\in I$,
\[
\gamma'(t)=\mathbf{X}(\gamma(t)).
\]
If $0\in I$, we call $\gamma(0)$ the \textit{initial point of $\gamma$}.
\end{definition}
\begin{example}\label{ej:nociones-fundamentales-consideramos-coordenadas-estandar-curvas-integrales-precisamen}
 \begin{enumerate}
 \item Consider $\mathbb{R}^{2}$ with standard coordinates $(x,y)$, and let $\mathbf{X}=\displaystyle\frac{\partial}{\partial x}$. The integral curves of $\mathbf{X}$ are precisely the lines parallel to the $x$-axis, parametrized as
 \[
 \gamma(t)=(a+t,b).
 \]
 In particular, there is a unique integral curve starting at each point of the plane, and the images of distinct integral curves either coincide or are disjoint.

 \item Let $\mathbf{W}=x\displaystyle\frac{\partial }{\partial y}-y\displaystyle\frac{\partial}{\partial x}$ on $\mathbb{R}^{2}$. If $\gamma\colon \mathbb{R}\longrightarrow \mathbb{R}^{2}$ is a smooth curve given by $\gamma(t)=(x(t),y(t))$, the condition $\gamma'(t)=\mathbf{W}(\gamma(t))$ becomes
 \[
 x'(t)\frac{\partial}{\partial x}\Biggr|_{\gamma(t)}+y'(t)\frac{\partial}{\partial y}\Biggr|_{\gamma(t)}
 =
 x(t)\frac{\partial }{\partial y}\Biggr|_{\gamma(t)}-y(t)\frac{\partial}{\partial x}\Biggr|_{\gamma(t)}.
 \]
 Equating components gives the system of ordinary differential equations
 \[
 \begin{cases}
 x'(t)=-y(t),\\
 y'(t)=x(t).
 \end{cases}
 \]
 Its solutions are
 \[
 x(t)=a\cos(t)-b\sin(t),\qquad
 y(t)=a\sin(t)+b\cos(t),
 \qquad a,b\in \mathbb{R}.
 \]
 Thus the integral curves are given by
 \[
 \gamma(t)=(a\cos(t)-b\sin(t),\,a\sin(t)+b\cos(t)).
 \]
 When $a=b=0$, we obtain the constant solution $\gamma(t)\equiv 0$; otherwise, the trajectories are circles centered at the origin, traversed counterclockwise. Since $\gamma(0)=(a,b)$, there is a unique integral curve starting at each point of $\mathbb{R}^{2}$, and the images of integral curves either coincide or are disjoint.
 \begin{center}
 \includegraphics[scale=0.9]{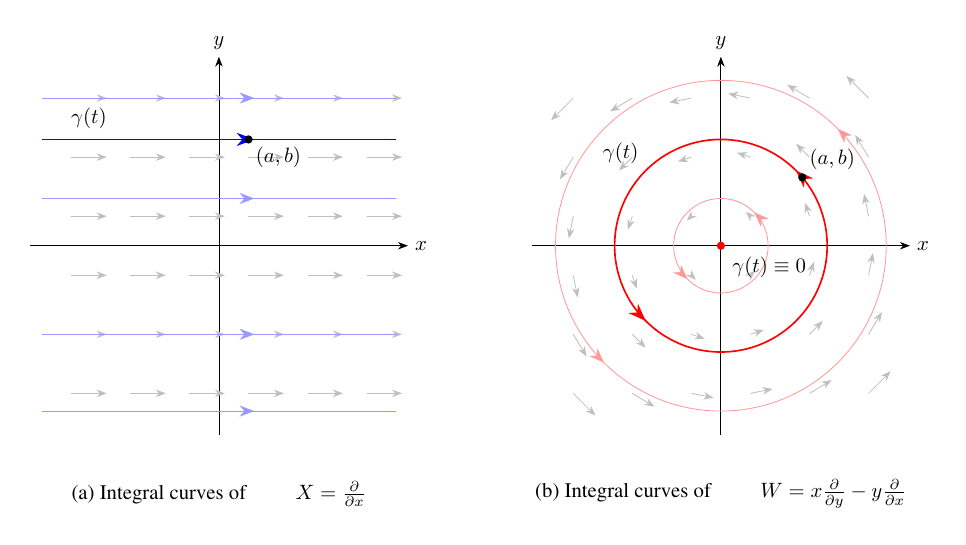}
 \end{center}
 \end{enumerate}
\end{example}
As the examples illustrate, finding integral curves ``reduces'' to solving a system of ordinary differential equations in a chart. Suppose that $\mathbf{X}\in \mathfrak{X}(M)$ and that $\gamma\colon I\longrightarrow M$ is a curve.

Let $(U,\phi)$ be a chart such that $\gamma(I)\subset U$, and write $\phi\circ \gamma=\widetilde{\gamma}=(x^{1}(t),\dots,x^{n}(t))$. Also write the vector field $\mathbf{X}$ in coordinates as $\mathbf{X} = X^{i}\boldsymbol{\partial}_{i}$.

The condition that $\gamma$ be an integral curve of $\mathbf{X}$, namely $\gamma'(t)=\mathbf{X}(\gamma(t))$, becomes in coordinates
\[
\dot{x}^{i}(t)\boldsymbol{\partial}_{i}\Biggr|_{\gamma(t)}
=
X^{i}(x(t))\boldsymbol{\partial}_{i}\Biggr|_{\gamma(t)}.
\]
Equating components gives the system of ordinary differential equations
\[
\begin{aligned}
\dot{x}^{1}(t)&=X^{1}(x^{1}(t),\dots,x^{n}(t)),\\
&\;\;\vdots\\
\dot{x}^{n}(t)&=X^{n}(x^{1}(t),\dots,x^{n}(t)).
\end{aligned}
\]

This shows that finding integral curves is equivalent, in coordinates, to solving an autonomous system of ordinary differential equations. The terminology comes precisely from this observation: finding integral curves corresponds to ``integrating'' the system.

The following result on ordinary differential equations will be fundamental for the subsequent development.
\begin{theorem}[Fundamental theorem for autonomous systems]\label{teo: fundamental sistemas autonomos}\index{fundamental theorem for autonomous systems@fundamental theorem for autonomous systems}
Let $U\subseteq \mathbb{R}^{n}$ be open and let $V\in C^{\infty}(U,\mathbb{R}^{n})$.
For $(t_{0},x_{0})\in \mathbb{R}\times U$, consider the initial value problem
\begin{equation}\label{eq:PVI-autonomo}
\begin{cases}
\dot y^{i}(t)=V^{i}(y^{1}(t),\dots,y^{n}(t)),\quad i\in\{1,\dots,n\},\\
y(t_{0})=x_{0}.
\end{cases}
\end{equation}

Then the following assertions hold:

\begin{enumerate}[label=(\alph*)]
\item (Local existence) For each $(t_{0},x_{0})\in \mathbb{R}\times U$, there exist an open interval $J_{0}\subseteq \mathbb{R}$ with $t_{0}\in J_{0}$ and an open set $U_{0}\subseteq U$ with $x_{0}\in U_{0}$ such that, for each $x\in U_{0}$, problem \eqref{eq:PVI-autonomo} with initial condition $y(t_{0})=x$ admits a solution $y\colon J_{0}\longrightarrow U$ of class $C^{1}$.

\item (Uniqueness) If $y_{1}$ and $y_{2}$ are two solutions of \eqref{eq:PVI-autonomo} defined on intervals containing $t_{0}$ and having the same initial condition, they agree on the intersection of their domains.

\item (Smooth dependence) If $J_{0}$ and $U_{0}$ are as in (a), the map
\[
\theta\colon J_{0}\times U_{0}\longrightarrow U,\qquad \theta(t,x)=y(t),
\]
where $y$ is the unique solution of \eqref{eq:PVI-autonomo} with $y(t_{0})=x$, is of class $C^{\infty}$.
\end{enumerate}
\end{theorem}
\begin{proof}
Apply Theorem~\ref{teo:picard-lindelof-banach} with
\(E=\mathbb R^n\) and \(F(t,y)=V(y)\). Since \(V\) is of class \(C^1\), it is
locally Lipschitz; the theorem gives a unique local solution.
Proposition~\ref{prop:unicidad-soluciones-edo-banach} extends uniqueness
to every common interval.

Theorem~\ref{teo:dependencia-parametros-edo-banach}, applied without
additional parameters, states that the natural domain
\[
\mathcal E:=\{(t,a)\in\mathbb R\times U:
t\text{ belongs to the maximal interval of the solution with }y(t_0)=a\}
\]
is open and that \((t,a)\mapsto y(t;t_0,a)\) is smooth. Since
\((t_0,x_0)\in\mathcal E\), there exist an interval \(J_0\) and an open set \(U_0\), with \(t_0\in J_0\), \(x_0\in U_0\), and \(J_0\times U_0\subseteq\mathcal E\). This gives
the uniform choice in (a), uniqueness in (b), and smooth dependence in
(c).
\end{proof}
We now establish existence of integral curves.
\begin{proposition}\label{prop:nociones-fundamentales-curva-suave-tal-curva-integral-comienza}
 Let $M$ be a smooth manifold without boundary and let
 $\mathbf{X}\in \mathfrak{X}(M)$. For each $p\in M$, there exist
 $\varepsilon>0$ and a smooth curve
 $\gamma\colon (-\varepsilon,\varepsilon)\longrightarrow M$ such that
 $\gamma$ is an integral curve of $\mathbf{X}$ starting at $p$.
\end{proposition}
\begin{proof}
Let $p\in M$. Take a smooth chart $(U,\phi)$ such that $p\in U$, and write $\phi=(x^{1},\dots,x^{n})$. Since $\mathbf{X}\in \mathfrak{X}(M)$, on $U$ it can be written as
$\mathbf{X}=X^{i}\boldsymbol{\partial}_{i}$, where each $X^{i}\colon \phi(U)\longrightarrow \mathbb{R}$ is smooth.

Consider the map $V\colon \phi(U)\longrightarrow \mathbb{R}^{n}$ given by
$V(x)=\bigl(X^{1}(x),\dots,X^{n}(x)\bigr)$, which is smooth. Let $x_{0}:=\phi(p)\in \phi(U)$. By Theorem~\ref{teo: fundamental sistemas autonomos}, applied to the system
$\dot y^{i}(t)=V^{i}(y^{1}(t),\dots,y^{n}(t))$, $i\in\{1,\dots,n\}$,
with initial condition $y(0)=x_{0}$, there exist $\varepsilon>0$ and a smooth curve
$y\colon (-\varepsilon,\varepsilon)\longrightarrow \phi(U)$ such that
\[
\dot y^{i}(t)=X^{i}(y^{1}(t),\dots,y^{n}(t)),\quad i\in\{1,\dots,n\},\qquad y(0)=x_{0}.
\]

Define $\gamma\colon (-\varepsilon,\varepsilon)\longrightarrow M$ by $\gamma(t):=\phi^{-1}(y(t))$. Then $\gamma$ is smooth and $\gamma(0)=p$.

To verify that $\gamma$ is an integral curve of $\mathbf{X}$, note that $\phi\circ\gamma=y=(y^{1},\dots,y^{n})$, so
\[
\gamma'(t)
=
\dot y^{i}(t)\boldsymbol{\partial}_{i}\Biggr|_{\gamma(t)}
=
X^{i}(y(t))\boldsymbol{\partial}_{i}\Biggr|_{\gamma(t)}
=
\mathbf{X}(\gamma(t)).
\]
Consequently, $\gamma$ is an integral curve of $\mathbf{X}$ starting at $p$.
\end{proof}
We now examine how affine reparametrizations affect integral curves:
 \begin{lemma}[rescaling]\label{lema: reescalamiento para curvas integrales}\index{rescaling lemma!for integral curves}
Let $M$ be a smooth manifold without boundary, let
$\mathbf{X}\in \mathfrak{X}(M)$, let $I\subseteq \mathbb{R}$ be an
interval, and let $\gamma\colon I\longrightarrow M$ be an integral curve of
$\mathbf{X}$. For any $a\in \mathbb{R}$, the curve
$\widetilde{\gamma}\colon \widetilde{I}\longrightarrow M$ defined by
$\widetilde{\gamma}(t)=\gamma(at)$ is an integral curve of
$a\mathbf{X}$, where
$\widetilde{I}=\{t\in \mathbb{R}\mid at\in I\}$.
\end{lemma}

\begin{proof}
Let $t_{0}\in \widetilde{I}$ and let $f$ be a smooth function defined on a neighborhood of $\widetilde{\gamma}(t_{0})$. Then
\[
\widetilde{\gamma}'(t_{0})f
=
\frac{d}{dt}\Biggr|_{t=t_{0}}(f\circ \widetilde{\gamma})(t)
=
\frac{d}{dt}\Biggr|_{t=t_{0}}(f\circ \gamma)(at)
=
a\frac{d}{dt}\Biggr|_{t=at_{0}}(f\circ \gamma)(t).
\]
Since $\gamma$ is an integral curve of $\mathbf{X}$, we have
\[
\frac{d}{dt}\Biggr|_{t=at_{0}}(f\circ \gamma)(t)
=
\mathbf{X}_{\gamma(at_{0})}f
=
\mathbf{X}_{\widetilde{\gamma}(t_{0})}f.
\]
Therefore,
\[
\widetilde{\gamma}'(t_{0})f
=
a\,\mathbf{X}_{\widetilde{\gamma}(t_{0})}f,
\]
proving that $\widetilde{\gamma}$ is an integral curve of $a\mathbf{X}$.
\end{proof}
\begin{lemma}[translation]\label{lema: traslacion curvas integrales}\index{translation lemma@translation lemma!for integral curves}
Let $M$ be a smooth manifold without boundary, let
$\mathbf{X}\in\mathfrak X(M)$, let $I\subseteq\mathbb R$ be an interval, and
let $\gamma\colon I\longrightarrow M$ be an integral curve of $\mathbf{X}$.
For any $b\in \mathbb{R}$, the curve
$\widetilde{\gamma}\colon \widetilde{I}\longrightarrow M$ given by
$\widetilde{\gamma}(t)=\gamma(t+b)$ is an integral curve of $\mathbf{X}$,
where $\widetilde{I}=\{t\in \mathbb{R}\mid t+b\in I\}$.
\end{lemma}

\begin{proof}
Let $t_{0}\in \widetilde{I}$ and let $f$ be a smooth function defined on a neighborhood of $\widetilde{\gamma}(t_{0})$. Then
\[
\widetilde{\gamma}'(t_{0})f
=
\frac{d}{dt}\Biggr|_{t=t_{0}}(f\circ \widetilde{\gamma})(t)
=
\frac{d}{dt}\Biggr|_{t=t_{0}}(f\circ \gamma)(t+b)
=
\frac{d}{dt}\Biggr|_{t=t_{0}+b}(f\circ \gamma)(t).
\]
Since $\gamma$ is an integral curve of $\mathbf{X}$, we have
\[
\frac{d}{dt}\Biggr|_{t=t_{0}+b}(f\circ \gamma)(t)
=
\mathbf{X}_{\gamma(t_{0}+b)}f
=
\mathbf{X}_{\widetilde{\gamma}(t_{0})}f.
\]
Therefore,
\[
\widetilde{\gamma}'(t_{0})f
=
\mathbf{X}_{\widetilde{\gamma}(t_{0})}f,
\]
proving that $\widetilde{\gamma}$ is an integral curve of $\mathbf{X}$.
\end{proof}
Smooth functions allow us to relate vector fields through the notion of $F$-related fields, which we introduce next.

\begin{definition}\label{def:nociones-fundamentales-variedad-7}\index{related vector fields}
Let $M$ and $N$ be smooth manifolds without boundary, let $F\colon M\longrightarrow N$ be a smooth map, and let $\mathbf{X}\in \mathfrak{X}(M)$ and $\mathbf{Y}\in \mathfrak{X}(N)$. We say that $\mathbf{X}$ and $\mathbf{Y}$ are \textit{$F$-related} if
\[
dF_p(\mathbf{X}_p)=\mathbf{Y}_{F(p)} \quad \text{for every } p\in M.
\]
\end{definition}
\begin{proposition}\label{prop: naturalidad de curvas integrales}
Let $M$ and $N$ be smooth manifolds without boundary and let $F\colon M\longrightarrow N$ be a smooth map. Then $\mathbf{X}\in\mathfrak{X}(M)$ and $\mathbf{Y}\in \mathfrak{X}(N)$ are $F$-related if and only if, for every integral curve $\gamma$ of $\mathbf{X}$, the curve $F\circ \gamma$ is an integral curve of $\mathbf{Y}$.
\end{proposition}

\begin{proof}
First suppose that $\mathbf{X}$ and $\mathbf{Y}$ are $F$-related, and let $\gamma\colon I\longrightarrow M$ be an integral curve of $\mathbf{X}$. Define $\sigma:=F\circ \gamma$. Then
\[
\sigma'(t)
=
(F\circ \gamma)'(t)
=
dF_{\gamma(t)}(\gamma'(t))
=
dF_{\gamma(t)}(\mathbf{X}_{\gamma(t)})
=
\mathbf{Y}_{F(\gamma(t))}
=
\mathbf{Y}_{\sigma(t)},
\]
so $\sigma$ is an integral curve of $\mathbf{Y}$.

Conversely, suppose $F$ sends integral curves of $\mathbf{X}$ to integral curves of $\mathbf{Y}$. Let $p\in M$ and let $\gamma\colon (-\varepsilon,\varepsilon)\longrightarrow M$ be an integral curve of $\mathbf{X}$ with $\gamma(0)=p$. Then $F\circ \gamma$ is an integral curve of $\mathbf{Y}$ with $(F\circ \gamma)(0)=F(p)$, and therefore
\[
\mathbf{Y}_{F(p)}
=
(F\circ \gamma)'(0)
=
dF_{p}(\gamma'(0))
=
dF_{p}(\mathbf{X}_{p}).
\]
This shows that $\mathbf{X}$ and $\mathbf{Y}$ are $F$-related.
\end{proof}
\begin{figura}[h]
    \includegraphics[scale=0.65]{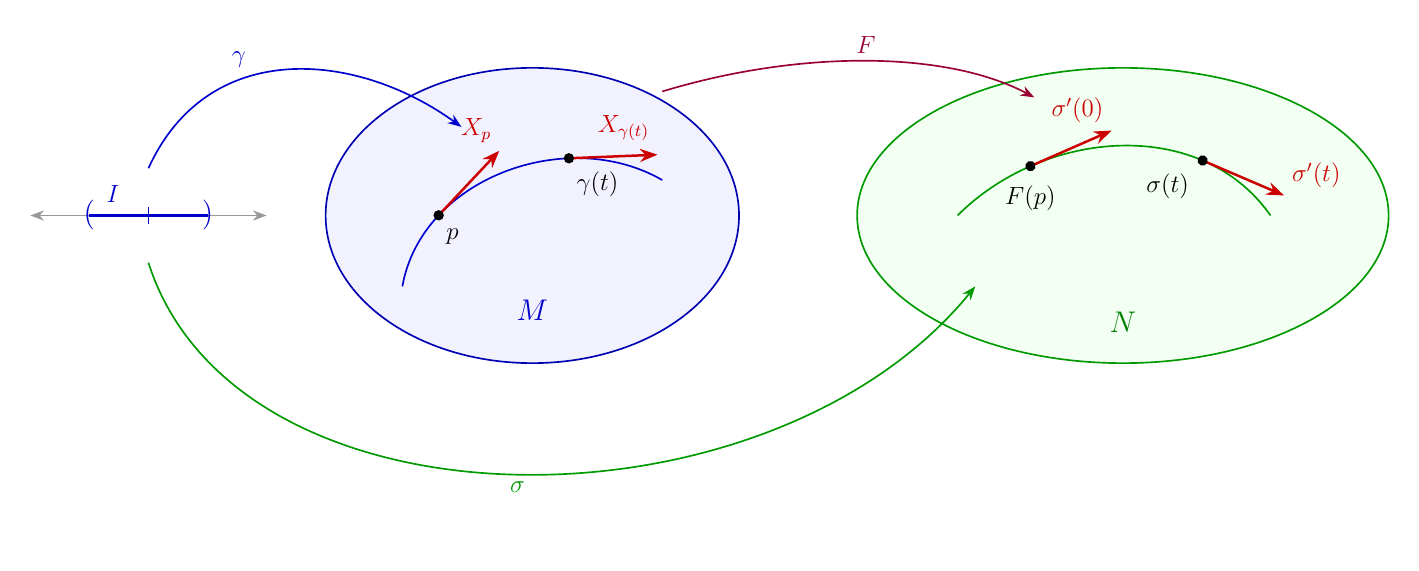}
    \caption{$F$-related vector fields.}
    \label{fig:campos-F-relacionados}
\end{figura}
We now study the notion of a \textit{flow}. To do so, we focus on the family of integral curves associated with a vector field.

Let $M$ be a smooth manifold without boundary and let $\mathbf{X}\in \mathfrak{X}(M)$. Suppose that, for each $p\in M$, $\mathbf{X}$ has a unique integral curve defined for every $t\in \mathbb{R}$ and starting at $p$, which we denote by $\theta^{(p)}\colon \mathbb{R}\longrightarrow M$.

For each $t\in \mathbb{R}$, define a map $\theta_{t}\colon M\longrightarrow M$ by $\theta_{t}(p):=\theta^{(p)}(t)$. Intuitively, each $\theta_t$ moves points of $M$ along their integral curves for time $t$.

The translation lemma implies that $t\mapsto \theta^{(p)}(t+s)$ is an integral curve of $\mathbf{X}$ starting at $q=\theta^{(p)}(s)$. Uniqueness of integral curves gives $\theta^{(q)}(t)=\theta^{(p)}(t+s)$. In terms of the maps $\theta_t$, this becomes
\[
\theta_t\circ \theta_s=\theta_{t+s}.
\]
Moreover, $\theta_0(p)=\theta^{(p)}(0)=p$, so the map $\theta\colon \mathbb{R}\times M\longrightarrow M$, given by $\theta(t,p)=\theta_t(p)$, defines an action of the additive group $\mathbb{R}$ on $M$.

With this in mind, we give the following definition.

\begin{definition}\label{def:nociones-fundamentales-flujo-global-en}\index{global flow}
Let $M$ be a smooth manifold without boundary. A \textit{\textbf{global flow on $M$}} is a continuous action of the additive group $\mathbb{R}$ on $M$, that is, a map $\theta\colon \mathbb{R}\times M\longrightarrow M$ such that, for any $s,t\in \mathbb{R}$ and $p\in M$,
\[
\theta(t,\theta(s,p))=\theta(t+s,p),\qquad \theta(0,p)=p.
\]
If $\theta$ is smooth, we say that $\theta$ is a \textit{\textbf{smooth global flow on $M$}}.

Given a global flow $\theta$ on $M$, define:
\begin{enumerate}
 \item For each $t\in \mathbb{R}$, the map $\theta_{t}\colon M\longrightarrow M$ given by $\theta_{t}(p)=\theta(t,p)$. Then $\theta_t\circ \theta_s=\theta_{t+s}$ and $\theta_0=\operatorname{Id}_M$, so each $\theta_t$ is a homeomorphism; if $\theta$ is smooth, $\theta_t$ is a diffeomorphism.

 \item For each $p\in M$, the curve $\theta^{(p)}\colon \mathbb{R}\longrightarrow M$ given by $\theta^{(p)}(t)=\theta(t,p)$. The image of this curve is called the orbit of $p$.
\end{enumerate}
\end{definition}
The following proposition shows that every smooth global flow arises from the integral curves of a smooth vector field.
\begin{definition}\label{def:nociones-fundamentales-generador-infinitesimal-de}\index{infinitesimal generator}
 Let $M$ be a smooth manifold without boundary and let
 $\theta\colon \mathbb{R}\times M\longrightarrow M$ be a smooth global
 flow. For each $p\in M$, define a tangent vector
 $\mathbf{X}_{p}\in T_{p}M$ by
 \[\mathbf{X}_{p}=(\theta^{(p)})'(0).\]
 The assignment $p\mapsto \mathbf{X}_{p}$ is a rough vector field on
 $M$, called the \textit{\textbf{infinitesimal generator of $\theta$}}.
\end{definition}
\begin{proposition}\label{prop:nociones-fundamentales-flujo-global-suave-definido-generador-infinitesimal}
 Let $M$ be a smooth manifold without boundary and let
 $\theta\colon \mathbb{R}\times M\longrightarrow M$ be a smooth global flow.
 The infinitesimal generator $\mathbf{X}$ of $\theta$ is a smooth vector
 field on $M$, and each $\theta^{(p)}$ is an integral curve of $\mathbf{X}$.
\end{proposition}
\begin{proof}
 To see that $\mathbf{X}$ is smooth, Proposition~\ref{prop: criterios de suavidad para un campo vectorial suave} says it suffices to show that $\mathbf{X}f$ is smooth for every $f\in C^{\infty}(U)$, with $U\subseteq M$ an arbitrary open set. We have \[\mathbf{X}f(p)=\mathbf{X}_{p}f=(\theta^{(p)})'(0)f=\frac{d}{dt}\Biggr|_{t=0}f(\theta^{(p)}(t))=\frac{\partial}{\partial t}\Biggr|_{(0,p)}f(\theta(t,p)).\]

 Since $f(\theta(t,p))$ is a smooth function of $(t,p)$ by composition, its partial derivative with respect to $t$ is also smooth. Thus $\mathbf{X}f(p)$ depends smoothly on $p$, so $\mathbf{X}$ is smooth. Finally, to see that the curve is an integral curve of $\mathbf{X}$, we must prove that $(\theta^{(p)})'(t)=\mathbf{X}_{\theta^{(p)}(t)}$ for every $p\in M$ and every $t\in \mathbb{R}$. Let $t_{0}\in \mathbb{R}$ be arbitrary and let $q=\theta^{(p)}(t_{0})=\theta_{t_{0}}(p)$. The group laws satisfied by $\theta$ give, for every $t\in \mathbb{R}$,
 \[\theta^{(q)}(t)=\theta_{t}(q)=\theta_{t}(\theta_{t_{0}}(p))=\theta_{t+t_{0}}(p)=\theta^{(p)}(t+t_{0}).\]
Therefore, for any $f\in C^{\infty}(U)$ where $q\in U$, \[\mathbf{X}_{q}f=(\theta^{(q)})'(0)f=\frac{d}{dt}\Biggr|_{t=0}f(\theta^{(q)}(t))=\frac{d}{dt}\Biggr|_{t=0}f(\theta^{(p)}(t+t_{0}))=(\theta^{(p)})'(t_{0})f.\]

\end{proof}
We have seen that every smooth global flow gives rise to a smooth vector field whose integral curves are exactly the curves defined by the flow. Conversely, we would like to say that every smooth vector field is the infinitesimal generator of a smooth global flow. This is not always true, however, since some smooth vector fields have integral curves that are not defined for every $t\in \mathbb{R}$, as the following example shows:
\begin{example}\label{ej:nociones-fundamentales-campo-vectorial-curva-integral-comienza-esta}
 Let $M=\mathbb{R}^{2}$ and consider the vector field
 \[\mathbf{X}(x,y)=x^{2}\frac{\partial}{\partial x}.\]
 The integral curve starting at $(1,0)$ is given by
 \[\gamma(t)=\left(\frac{1}{1-t},0\right),\]
 which solves the differential equation $x'=x^{2}$ with initial condition $x(0)=1$.

 Note that $\gamma(t)$ cannot be extended beyond $t=1$, since $\displaystyle\frac{1}{1-t}\to +\infty$ as $t\to 1^{-}$. Thus the integral curve is not defined for every $t\in\mathbb{R}$, and consequently $\mathbf{X}$ does not generate a global flow.
\end{example}
For this reason, we need to restrict the domains on which flows are defined while retaining their characteristic properties.
\begin{definition}\label{def:nociones-fundamentales-flujo}\index{flow}
 Let $M$ be a smooth manifold without boundary. A flow domain for $M$ is an open subset $\mathcal{D}\subseteq \mathbb{R}\times M$ such that, for each $p\in M$, the set \[\mathcal{D}^{(p)}=\{t\in \mathbb{R}\mid (t,p)\in \mathcal{D}\}\] is an open interval containing $0$. A \textit{\textbf{flow}} on $M$ is a continuous function $\theta\colon \mathcal{D}\longrightarrow M$, where $\mathcal{D}\subseteq \mathbb{R}\times M$ is a flow domain, satisfying the following group laws:
 \begin{enumerate}
 \item For every $p\in M$, $\theta(0,p)=p$.
 \item For every $s\in \mathcal{D}^{(p)}$ and every $t\in \mathcal{D}^{(\theta(s,p))}$ such that $s+t\in \mathcal{D}^{(p)}$, $\theta(t,\theta(s,p))=\theta(t+s,p)$.
 \end{enumerate}
 We sometimes call $\theta$ a \textit{\textbf{local flow}} to distinguish it from global flows. If $\theta$ is a flow, define $\theta_{t}(p)=\theta^{(p)}(t)=\theta(t,p)$ if $(t,p)\in \mathcal{D}$. For $t\in \mathbb{R}$, define \[M_{t}:=\{p\in M\mid (t,p)\in \mathcal{D}\},\] so that $p\in M_{t}\iff t\in \mathcal{D}^{(p)}\iff (t,p)\in \mathcal{D}$. If $\theta$ is smooth, the infinitesimal generator of $\theta$ is defined by $\mathbf{X}_{p}:=(\theta^{(p)})'(0)$.
\end{definition}
\begin{proposition}\label{prop:nociones-fundamentales-flujo-suave-generador-infinitesimal-campo-vectorial}
 Let $M$ be a smooth manifold without boundary. If
 $\theta\colon \mathcal{D}\longrightarrow M$ is a smooth flow, the
 infinitesimal generator $\mathbf{X}$ of $\theta$ is a smooth vector
 field, and each curve $\theta^{(p)}$ is an integral curve of $\mathbf{X}$.
\end{proposition}
\begin{proof}
Since $\mathcal D$ is open and contains $\{0\}\times M$, the map
\[
p\longmapsto \mathbf{X}_p:=d\theta_{(0,p)}(1,0)
\]
is a smooth section of $TM$: in a chart of $M$, its components are the
partial derivatives with respect to time of the components of
$\theta$, evaluated at $(0,p)$. This proves smoothness of $\mathbf{X}$.

Fix $p\in M$ and $t_0\in\mathcal D^{(p)}$. Since
$\mathcal D^{(\theta(t_0,p))}$ is an open interval containing zero,
for $s$ sufficiently close to zero the flow law gives
\[
\theta\bigl(s,\theta(t_0,p)\bigr)=\theta(t_0+s,p).
\]
Differentiating with respect to $s$ at $s=0$, we obtain
\[
\mathbf{X}_{\theta(t_0,p)}
=
\left.\frac{d}{ds}\right|_{s=0}\theta(t_0+s,p)
=
\frac{d}{dt}\theta^{(p)}(t_0).
\]
Since $t_0$ was arbitrary, $\theta^{(p)}$ satisfies the integral-curve
equation throughout its domain $\mathcal D^{(p)}$.
\end{proof}
Existence and uniqueness of the flow induced by a vector field is the most important result in flow theory. It will be useful later when we consider certain flows in the study of Sobolev spaces.
\begin{theorem}[fundamental theorem of flows]\label{teo: fundamental de flujos}\index{fundamental theorem of flows}
Let $M$ be a smooth manifold without boundary and let
$\mathbf{X}\in \mathfrak{X}(M)$. Then there exists a unique maximal smooth
flow $\theta\colon \mathcal{D}\longrightarrow M$ whose infinitesimal
 generator is $\mathbf{X}$ and which has the following properties:
\begin{enumerate}
 \item For each $p\in M$, the curve $\theta^{(p)}\colon \mathcal{D}^{(p)}\longrightarrow M$ is the unique maximal integral curve of $\mathbf{X}$ starting at $p$.
 \item If $s\in \mathcal{D}^{(p)}$, then $\mathcal{D}^{(\theta(s,p))}$ is the interval $\mathcal{D}^{(p)}-s=\{t-s\mid t\in \mathcal{D}^{(p)}\}$.
 \item For each $t\in \mathbb{R}$, the set $M_{t}$ is open in $M$, and $\theta_{t}\colon M_{t}\longrightarrow M_{-t}$ is a diffeomorphism with inverse $\theta_{-t}$.
\end{enumerate}
\end{theorem}
\begin{proof}
For each \(p\in M\), local existence and uniqueness in charts allow us
to join all integral curves starting at \(p\). This gives a
unique maximal integral curve
\[
\gamma_p\colon I_p\longrightarrow M,\qquad
0\in I_p,\qquad \gamma_p(0)=p.
\]
Define
\[
\mathcal D:=\{(t,p)\in\mathbb R\times M:t\in I_p\},
\qquad
\theta(t,p):=\gamma_p(t).
\]

Fix \((t_1,p_1)\in\mathcal D\). The image of the compact portion of
\(I_{p_1}\) between \(0\) and \(t_1\) is covered by finitely many
chart domains. In each chart,
Theorem~\ref{teo: fundamental sistemas autonomos} provides a common existence
interval and smooth dependence on the initial data.
Divide the interval between \(0\) and \(t_1\) into finitely many subintervals
subordinate to those domains. On the first, the final value depends
smoothly on the initial point; this value serves as the initial data for the second,
and so on. Composing finitely many of these solution maps
gives an open set \(U_1\) and an interval \(J_1\), with \(p_1\in U_1\) and \(t_1\in J_1\), such that
\[
J_1\times U_1\subseteq\mathcal D
\]
such that \(\theta\) is smooth on \(J_1\times U_1\). Therefore,
\(\mathcal D\) is open and \(\theta\) is smooth.

Let \(s\in I_p\) and \(q=\gamma_p(s)\). The curves
\[
t\longmapsto\gamma_p(t+s),
\qquad
t\longmapsto\gamma_q(t)
\]
are integral curves of \(\mathbf X\) and take the value \(q\) at \(t=0\). By
uniqueness, they agree on their common domain. Maximality in both
directions implies
\[
I_q=I_p-s,
\qquad
\gamma_q(t)=\gamma_p(t+s),\quad t\in I_p-s.
\]
This proves (2), the flow laws, and
\(\theta^{(p)}=\gamma_p\), which is (1). Moreover,
\[
\left.\frac d{dt}\right|_{t=0}\theta(t,p)
=\gamma_p'(0)=\mathbf X_p,
\]
so the generator is \(\mathbf X\).

Since \(\mathcal D\) is open, each
\(M_t=\{p:(t,p)\in\mathcal D\}\) is open. If \(p\in M_t\) and
\(q=\theta_t(p)\), then \(I_q=I_p-t\), so \(q\in M_{-t}\), and
the flow law gives \(\theta_{-t}(q)=p\). The reverse argument starts
with \(q\in M_{-t}\). Thus
\(\theta_t\colon M_t\to M_{-t}\) is a diffeomorphism with inverse
\(\theta_{-t}\), proving (3).

Any other smooth flow with generator \(\mathbf X\) has integral curves
of \(\mathbf X\) as its orbits. If it is maximal, each orbit agrees
with \(\gamma_p\); hence both the domain and the flow agree. This
proves uniqueness.
\end{proof}
\begin{figura}[h]
    \includegraphics[scale=0.9]{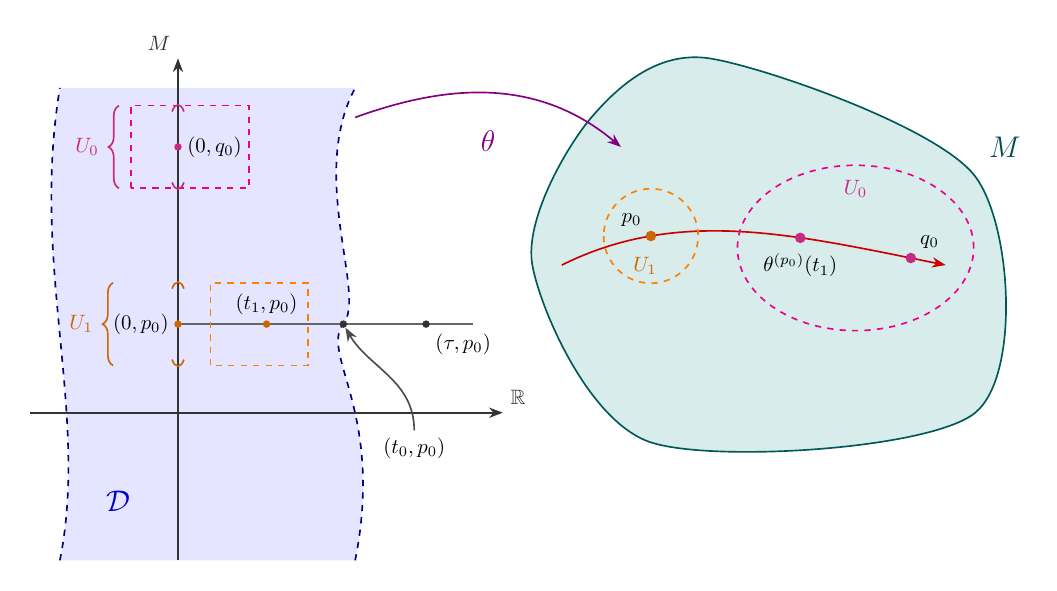}
    \caption{Argument showing that $\mathcal{D}$ is open.}
    \label{fig:D-abierto}
\end{figura}The flow whose existence is guaranteed by Theorem~\ref{teo: fundamental de flujos} is called the \textit{flow generated by $\mathbf{X}$}. The term \textit{infinitesimal generator} is justified by the fact that the vector field $\mathbf{X}$ describes the instantaneous variation of the flow, in the sense that
\[
\mathbf{X}(p)=\left.\frac{d}{dt}\theta(t,p)\right|_{t=0}.
\]
In a smooth chart, this means that, at each point, the vector $\mathbf{X}(p)$ determines the tangent direction of the integral curve through $p$, and the flow is obtained by integrating these directions over time. Flows are natural with respect to smooth functions in the following sense:
\begin{proposition}\label{prop:nociones-fundamentales-variedades-suaves-mapeo-suave-flujo-flujo}
 Let $M$ and $N$ be smooth manifolds without boundary, let $F\colon M\longrightarrow N$ be a smooth map, and let $\mathbf{X}\in \mathfrak{X}(M)$ and $\mathbf{Y}\in \mathfrak{X}(N)$. Let $\theta$ be the flow of $\mathbf{X}$ and $\eta$ the flow of $\mathbf{Y}$. If $\mathbf{X}$ and $\mathbf{Y}$ are $F$-related, then, for each $t\in \mathbb{R}$, $F(M_{t})\subseteq N_{t}$ and $\eta_{t}\circ F = F\circ \theta_{t}$ on $M_{t}$: \[
\begin{tikzcd}
M_t \arrow[r, "\theta_t"] \arrow[d, "F"'] & M_{-t} \arrow[d, "F"] \\
N_t \arrow[r, "\eta_t"] & N_{-t}
\end{tikzcd}
\]
\end{proposition}
\begin{proof}
 By Proposition~\ref{prop: naturalidad de curvas integrales}, for any $p\in M$, the curve $F\circ\theta^{(p)}$ is an integral curve of
 $\mathbf{Y}$ starting at $F(\theta^{(p)}(0))=F(p)$. By uniqueness of integral curves, the maximal integral curve $\eta^{(F(p))}$ must be defined at least on $\mathcal{D}^{(p)}$, and $F\circ \theta^{(p)}=\eta^{(F(p))}$ on that interval. This means that \[p\in M_{t}\Rightarrow t\in \mathcal{D}^{(p)}\Rightarrow t\in \mathcal{D}^{(F(p))}\Rightarrow F(p)\in N_{t},\] equivalently $F(M_{t})\subseteq N_{t}$, and \[F(\theta^{(p)}(t))=\eta^{(F(p))}(t),\qquad \forall t\in \mathcal{D}^{(p)},\] equivalently $\eta_{t}(F(p))=F(\theta_{t}(p))$ for every $p\in M_{t}$.
\end{proof}
This gives the following corollary:
\begin{corollary}[Invariance of flows under diffeomorphisms]\label{cor: invariancia de flujos bajo difeos}\index{invariance of flows under diffeomorphisms}
Let $M$ and $N$ be smooth manifolds without boundary, let
$F\colon M\longrightarrow N$ be a diffeomorphism, and let
$\mathbf{X}\in \mathfrak{X}(M)$ have flow $\theta$. Then the vector
field $F_{*}\mathbf{X}\in \mathfrak{X}(N)$ defined by
\[
(F_*\mathbf{X})(q)=dF_{F^{-1}(q)}\bigl(\mathbf{X}(F^{-1}(q))\bigr)
\]
has flow given by
\[
\eta_t = F\circ \theta_t \circ F^{-1},
\]
with domain $N_t = F(M_t)$ for each $t\in \mathbb{R}$.
\end{corollary}
\begin{proof}
Set \(\mathbf Y:=F_*\mathbf X\). By the definition of the pushforward field,
\[
dF_p(\mathbf X(p))
=
\mathbf Y(F(p)),
\qquad p\in M,
\]
so \(\mathbf X\) and \(\mathbf Y\) are \(F\)-related. Proposition~\ref{prop:nociones-fundamentales-variedades-suaves-mapeo-suave-flujo-flujo} applied to \(F\) gives
\[
F(M_t)\subseteq N_t,
\qquad
\eta_t\circ F=F\circ\theta_t
\quad\text{on }M_t.
\]
Apply the same proposition to \(F^{-1}\), observing that
\((F^{-1})_*\mathbf Y=\mathbf X\). This gives
\[
F^{-1}(N_t)\subseteq M_t.
\]
Applying \(F\) to this inclusion yields \(N_t\subseteq F(M_t)\); therefore,
\[
N_t=F(M_t).
\]
Finally, composing \(\eta_t\circ F=F\circ\theta_t\) on the right with \(F^{-1}\) gives
\[
\eta_t=F\circ\theta_t\circ F^{-1}
\quad\text{on }N_t.
\]
\end{proof}

 We record another important property of integral curves:
 \begin{lemma}[escape]\label{lema: de escape}\index{escape lemma}
Let $M$ be a smooth manifold without boundary and let $\mathbf{X}\in \mathfrak{X}(M)$. If
$\gamma\colon J\longrightarrow M$ is a maximal integral curve of $\mathbf{X}$ and its
right endpoint $b:=\displaystyle\sup J$ is finite, then, for every $t_{0}\in J$, the
tail $\gamma([t_{0},b))$ is not contained in any compact subset of $M$.
The assertion for the left endpoint follows by reversing the parameter:
if $a:=\displaystyle\inf J> -\infty$, no tail
$\gamma((a,t_{0}])$ is contained in a compact set.
 \end{lemma}
 \begin{proof}
Suppose that $\gamma([t_0,b))\subseteq K$ with $K$ compact. For a
sequence $t_j\to b^{-}$, compactness gives a subsequence, not
relabeled, such that $\gamma(t_j)\to p\in K$. The local existence theorem
and dependence on the initial data provide a neighborhood $U$ of $p$
and a common $\varepsilon>0$ for all initial data in $U$. For sufficiently large $j$,
$\gamma(t_j)\in U$ and $b-t_j<\tfrac{\varepsilon}{2}$. The local solution
starting at $\gamma(t_j)$ at time $t_j$ agrees by uniqueness with
$\gamma$ on their common domain and extends it to $t_j+\varepsilon>b$,
contradicting maximality. For the left endpoint, apply the same
argument to the field $-\mathbf{X}$ and the time-reversed curve.
 \end{proof}
 We give some local estimates for partial derivatives of a flow, which will be useful later when studying Sobolev spaces on manifolds of bounded geometry.
 \begin{lemma}[Grönwall's inequality]\label{lema: gronwall}\index{Gronwall inequality@Grönwall's inequality}
Let $u\colon [0,T]\longrightarrow [0,\infty)$ be a continuous function, let $a\geq 0$, and let $b\colon [0,T]\longrightarrow [0,\infty)$ be a continuous function. Suppose that
\[
u(t)\leq a+\int_{0}^{t} b(\tau)\,u(\tau)\,d\tau,
\qquad \forall t\in [0,T].
\]
Then
\[
u(t)\leq a\,\exp\left(\int_{0}^{t} b(\tau)\,d\tau\right),
\qquad \forall t\in [0,T].
\]
\end{lemma}
\begin{proof}
Define
\[
v(t):=a+\int_0^t b(\tau)u(\tau)\,d\tau.
\]
Then \(u\leq v\), \(v(0)=a\), and
\[
v'(t)=b(t)u(t)\leq b(t)v(t).
\]
If \(a>0\), we have \(v(t)>0\), so
\[
\frac{d}{dt}\log v(t)\leq b(t).
\]
Integrating between \(0\) and \(t\),
\[
v(t)\leq a\exp\left(\int_0^tb(\tau)\,d\tau\right),
\]
and the same bound holds for \(u\). If \(a=0\), apply the preceding
argument, for each \(\varepsilon>0\), to
\[
v_\varepsilon(t):=\varepsilon+\int_0^tb(\tau)u(\tau)\,d\tau.
\]
This gives
\[
u(t)\leq v_\varepsilon(t)
\leq\varepsilon\exp\left(\int_0^tb(\tau)\,d\tau\right).
\]
Letting \(\varepsilon\to 0^{+}\) yields \(u(t)=0\), which is the
assertion in this case.
\end{proof}
\begin{lemma}\label{lema: estimacion derivadas flujo variedad}
Let $M$ be a smooth manifold without boundary of dimension $n$, and let $\mathbf{X}\in \mathfrak{X}(M)$ be a smooth vector field. Let
\[
\Phi\colon D\subseteq \mathbb{R}\times M\longrightarrow M
\]
be the flow of $\mathbf{X}$.

Let $(U,\phi)$ be a chart of $M$ such that $\overline{U}$ is compact. Suppose that $p\in U$, that $t\geq0$, and that
\[
(\tau,p)\in D\quad\text{and}\quad \Phi(\tau,p)\in U,
\qquad \tau\in [0,t].
\]
Then, for each multi-index $\alpha\in \mathbb{N}_{0}^{n}$ with $|\alpha|\geq 1$, there exists a continuous function
\[
\operatorname{Expr}_{\alpha}\colon [0,\infty)\times [0,\infty)\longrightarrow [0,\infty),
\]
nondecreasing in each variable and depending only on $\alpha$ and $n$, such that
\[
\bigl\|\partial_{x}^{\alpha}(\phi\circ \Phi(t,\cdot)\circ \phi^{-1})(\phi(p))\bigr\|
\leq
\operatorname{Expr}_{\alpha}\left(
\displaystyle\sup_{\substack{\tau\in[0,t],\ \beta\in\mathbb N_0^n\\1\leq|\beta|\leq|\alpha|}}
\bigl\|\partial_{x}^{\beta}(\phi_*\mathbf{X})(\phi(\Phi(\tau,p)))\bigr\|,
\, t
\right).
\]

The vector norms are Euclidean. Moreover, if $|\alpha|=1$, we may take
\[
\operatorname{Expr}_{\alpha}(M,t)=e^{\sqrt n\,Mt},
\]
while for $|\alpha|\geq 2$ we can choose $\operatorname{Expr}_{\alpha}$ recursively.
\end{lemma}
\begin{proof}
Define
\[
V:=\phi(U)\subseteq \mathbb{R}^{n}
\]
and consider the vector field in coordinates
\[
\widetilde X:=\phi_{*}\mathbf{X}\in \mathfrak{X}(V).
\]
Since $\mathbf{X}$ is smooth and $\phi$ is a diffeomorphism onto its image, the field $\widetilde X$ is smooth.

Also define
\[
\widetilde{D}:=\bigl\{(\tau,x)\in \mathbb{R}\times V \bigm|
(\tau,\phi^{-1}(x))\in D\ \text{and}\
\Phi(\tau,\phi^{-1}(x))\in U\bigr\}
\]
and
\[
\widetilde{\Phi}\colon \widetilde{D}\longrightarrow V,
\qquad
\widetilde{\Phi}(\tau,x):=\phi\bigl(\Phi(\tau,\phi^{-1}(x))\bigr).
\]
The function $\widetilde{\Phi}$ is smooth. We will show that, for each $x\in V$, the curve
\[
\tau\longmapsto \widetilde{\Phi}(\tau,x)
\]
satisfies the differential equation associated with $\widetilde X$.

Fix $x\in V$ and write $q:=\phi^{-1}(x)$. Then
\[
\frac{\partial}{\partial \tau}\widetilde{\Phi}(\tau,x)
=
d\phi_{\Phi(\tau,q)}
\left(
\frac{\partial}{\partial \tau}\Phi(\tau,q)
\right).
\]
Since $\Phi$ is the flow of $\mathbf{X}$,
\[
\frac{\partial}{\partial \tau}\Phi(\tau,q)=\mathbf{X}_{\Phi(\tau,q)}.
\]
Therefore,
\[
\frac{\partial}{\partial \tau}\widetilde{\Phi}(\tau,x)
=
d\phi_{\Phi(\tau,q)}(\mathbf{X}_{\Phi(\tau,q)})
=
\widetilde X(\widetilde{\Phi}(\tau,x)).
\]
Moreover,
\[
\widetilde{\Phi}(0,x)=x.
\]
Thus $\tau\mapsto \widetilde{\Phi}(\tau,x)$ solves
\[
z'(\tau)=\widetilde X(z(\tau)),
\qquad z(0)=x.
\]

Take
\[
x_0:=\phi(p), \qquad Z(\tau):=\widetilde{\Phi}(\tau,x_0).
\]
Then
\[
Z'(\tau)=\widetilde X(Z(\tau)),
\qquad Z(0)=x_0,
\]
and integration gives
\[
Z(t)=x_0+\int_0^t \widetilde X(Z(\tau))\,d\tau.
\]

\medskip

We proceed by induction on $|\alpha|$.

\medskip

\textbf{Case $|\alpha|=1$.}
Define
\[
Y_k(t):=\partial_{x_{0,k}}\widetilde{\Phi}(t,x_0).
\]
Differentiating the integral equation,
\[
Y_k(t)=e_k+\int_0^t d\widetilde X_{Z(\tau)}(Y_k(\tau))\,d\tau,
\]
and therefore
\[
Y_k'(t)=d\widetilde X_{Z(t)}(Y_k(t)), \qquad Y_k(0)=e_k.
\]

Let
\[
M_1(t):=\displaystyle\sup_{\substack{\tau\in[0,t]\\j\in\{1,\dots,n\}}}
\|\partial_j\widetilde X(Z(\tau))\|.
\]
For $z\in V$ and $v=(v_1,\dots,v_n)\in\mathbb R^n$, we have
\[
\|d\widetilde X_z(v)\|
\leq\max_{1\leq j\leq n}\|\partial_j\widetilde X(z)\|
\displaystyle\sum_{j=1}^n|v_j|
\leq\sqrt n\max_{1\leq j\leq n}\|\partial_j\widetilde X(z)\|\,\|v\|.
\]
Thus the operator norm of $d\widetilde X_{Z(\tau)}$ is at
most $\sqrt n\,M_1(t)$ for $0\leq\tau\leq t$. Fixing $t$,
\[
\|Y_k(s)\|\le 1+\int_0^s\sqrt n\,M_1(t)\|Y_k(\tau)\|\,d\tau,
\qquad 0\leq s\leq t.
\]
Applying Lemma~\ref{lema: gronwall} and evaluating at $s=t$,
\[
\|Y_k(t)\|\le e^{\sqrt n\,M_1(t)t}.
\]

\medskip

\textbf{Inductive step.}
Suppose $|\alpha|=m\ge 2$ and that the result holds for every multi-index of strictly smaller order.

Apply $\partial_{x_0}^{\alpha}$ to
\[
\frac{\partial}{\partial t}\widetilde{\Phi}(t,x_0)
=
\widetilde X(\widetilde{\Phi}(t,x_0)).
\]
Commuting derivatives,
\[
\frac{\partial}{\partial t}\bigl(\partial_{x_0}^{\alpha}\widetilde{\Phi}(t,x_0)\bigr)
=
\partial_{x_0}^{\alpha}(\widetilde X(\widetilde{\Phi}(t,x_0))).
\]

Applying the Faà di Bruno formula (Theorem~\ref{faa di bruno multivariable}), we obtain
\[
\partial_{x_0}^{\alpha}(\widetilde X(\widetilde{\Phi}(t,x_0)))
=
\displaystyle\sum_{\substack{\beta\in \mathbb{N}_{0}^{n}\\ 1\leq |\beta|\leq |\alpha|}}
\;
\displaystyle\sum_{\substack{
\gamma=(\gamma_{1},\dots,\gamma_{|\beta|})\\
|\gamma_i|>0,\ \displaystyle\sum_{i=1}^{|\beta|}\gamma_i=\alpha}}
\;
\displaystyle\sum_{l\in\{1,\dots,n\}^{|\beta|}}
c_{\alpha,\beta,\gamma,l}
\,
\partial_x^{\beta}\widetilde X(Z(t))
\prod_{i=1}^{|\beta|}
\partial_{x_0}^{\gamma_i}\widetilde{\Phi}_{l_i}(t,x_0).
\]

If $|\beta|=1$, then $\gamma=(\gamma_1)$ has only one component and $\gamma_1=\alpha$, giving
\[
d\widetilde X_{Z(t)}(\partial_{x_0}^{\alpha}\widetilde{\Phi}(t,x_0)).
\]

We therefore define explicitly
\[
G_\alpha(t):=
\displaystyle\sum_{\substack{\beta\in \mathbb{N}_{0}^{n}\\ 2\leq |\beta|\leq |\alpha|}}
\;
\displaystyle\sum_{\substack{
\gamma=(\gamma_{1},\dots,\gamma_{|\beta|})\\
|\gamma_i|>0,\ \displaystyle\sum_{i=1}^{|\beta|}\gamma_i=\alpha}}
\;
\displaystyle\sum_{l\in\{1,\dots,n\}^{|\beta|}}
c_{\alpha,\beta,\gamma,l}
\,
\partial_x^{\beta}\widetilde X(Z(t))
\prod_{i=1}^{|\beta|}
\partial_{x_0}^{\gamma_i}\widetilde{\Phi}_{l_i}(t,x_0).
\]

Thus
\[
\frac{\partial}{\partial t}(\partial_{x_0}^{\alpha}\widetilde{\Phi})
=
d\widetilde X_{Z(t)}(\partial_{x_0}^{\alpha}\widetilde{\Phi})
+
G_\alpha(t),
\qquad
\partial_{x_0}^{\alpha}\widetilde{\Phi}(0,x_0)=0.
\]

Define
\[
M_\alpha(t):=
\displaystyle\sup_{\substack{\tau\in[0,t],\ \beta\in\mathbb N_0^n\\1\leq|\beta|\leq|\alpha|}}
\|\partial_x^\beta \widetilde X(Z(\tau))\|.
\]

By the induction hypothesis,
\[
\|\partial_{x_0}^{\gamma}\widetilde{\Phi}(\tau,x_0)\|
\le \operatorname{Expr}_\gamma(M_\alpha(t),t).
\]

Taking norms in $G_\alpha$ and using these bounds shows that there exists a polynomial $P_\alpha$ with nonnegative coefficients such that
\[
\|G_\alpha(\tau)\|
\le
P_\alpha\bigl(M_\alpha(t),\operatorname{Expr}_{\gamma^{(1)}}(M_\alpha(t),t),\dots\bigr).
\]

Define
\[
B_\alpha(M,t):=P_\alpha\bigl(M,\operatorname{Expr}_{\gamma^{(1)}}(M,t),\dots\bigr).
\]

Let
\[
u(s):=\|\partial_{x_0}^{\alpha}\widetilde{\Phi}(s,x_0)\|.
\]
The norm comparison in the first-order case gives
$\|d\widetilde X_{Z(\tau)}\|_{\mathrm{op}}\leq\sqrt n\,M_\alpha(t)$
for $0\leq\tau\leq t$. Integrating the equation for
$\partial_{x_0}^{\alpha}\widetilde\Phi$, whose initial value is zero, gives
\[
u(s)\le tB_\alpha(M_\alpha(t),t)+\int_0^s\sqrt n\,M_\alpha(t)u(\tau)d\tau,
\qquad 0\leq s\leq t.
\]

Applying Lemma~\ref{lema: gronwall},
\[
u(s)\le tB_\alpha(M_\alpha(t),t)e^{\sqrt n\,M_\alpha(t)s}.
\]
Evaluating at $s=t$,
\[
\|\partial_{x_0}^{\alpha}\widetilde{\Phi}(t,x_0)\|
\le
tB_\alpha(M_\alpha(t),t)e^{\sqrt n\,M_\alpha(t)t}.
\]

Define
\[
\operatorname{Expr}_\alpha(M,t):=tB_\alpha(M,t)e^{\sqrt n\,Mt}.
\]

Finally, the function $\operatorname{Expr}_\alpha$ is nondecreasing in each variable. By the induction hypothesis, the functions $\operatorname{Expr}_\gamma$ have this property; since $P_\alpha$ has nonnegative coefficients, $B_\alpha$ is nondecreasing in each variable; and the function $(M,t)\mapsto t e^{\sqrt n\,Mt}$ is nondecreasing in $[0,\infty)\times[0,\infty)$. The assertion follows.

Finally,
\[
\partial_x^\alpha(\phi\circ \Phi(t,\cdot)\circ \phi^{-1})(\phi(p))
=
\partial_{x_0}^{\alpha}\widetilde{\Phi}(t,x_0),
\]
and
\[
\partial_x^\beta \widetilde X(Z(\tau))
=
\partial_x^\beta(\phi_*\mathbf{X})(\phi(\Phi(\tau,p))),
\]
which completes the proof.
\end{proof}

\section{Lie groups and smooth actions}
\label{sec:grupos-lie-acciones-fundamentos}

Many geometric symmetries depend on continuous parameters.
Translations of $\mathbb R^n$, rotations of the circle, and changes of
frame can be composed and inverted, but they can also be differentiated
with respect to those parameters. Lie groups bring these two structures together:
the group law expresses composition of symmetries, while the smooth structure
allows us to study how they vary.

Basic examples include $(\mathbb R^n,+)$, the circle $S^1$,
$\operatorname{GL}(n,\mathbb R)$, and its orthogonal subgroups; in the complex
case, unitary groups arise in the same way. Isometry
groups provide further examples when the necessary geometric
structure is available. Near the identity, the tangent space $T_eG$ and its
bracket form the Lie algebra of $G$. This infinitesimal linearization
allows us to study curves of symmetries through equations in a vector
space.

The relation between the group and its algebra will appear in flows and
one-parameter subgroups. Smooth actions will later lead to homogeneous
spaces and principal bundles; in infinite dimensions, analogous ideas
will arise in gauge groups and diffeomorphism groups.

In this section, we specify the Lie group conventions used
later. If $G$ is a group and
$g\in G$, we write
\[
 L_g(h):=gh,\qquad R_{\mathbf{g}}(h):=hg,\qquad C_g(h):=ghg^{-1}.
\]
For the standard finite-dimensional theory and a comparison of
conventions, see \cite[Chs.~7 and~20]{LeeS}. We will explicitly
indicate external results whose proofs are not reproduced.

\begin{definition}[Lie group]
\label{def:grupo-lie-finito-dimensional}
\index{Lie group}
A \textit{finite-dimensional real Lie group} is a smooth manifold
$G$ without boundary, equipped with a group structure for which the
maps
\[
 \mu\colon G\times G\longrightarrow G,\qquad \mu(g,h)=gh,
 \qquad
 \iota\colon G\longrightarrow G,\qquad \iota(g)=g^{-1},
\]
are smooth. Denote the identity element of $G$ by $e$.
A \textit{Lie group homomorphism} is a group homomorphism that
is also smooth.
\end{definition}

\begin{proposition}
\label{prop:traslaciones-conjugaciones-grupos-lie}
Let $G$ and $H$ be finite-dimensional real Lie groups whose underlying smooth
manifolds have no boundary. For every $g\in G$, the maps
$L_g$, $R_{\mathbf{g}}$, and $C_g$ are diffeomorphisms, and
\[
 L_g^{-1}=L_{g^{-1}},\qquad
 R_{\mathbf{g}}^{-1}=R_{g^{-1}},\qquad
 C_g^{-1}=C_{g^{-1}}.
\]
Moreover,
\[
 C_{gh}=C_g\circ C_h.
\]
If $\Phi\colon G\to H$ is a Lie group homomorphism, then
\[
 \Phi\circ L_g=L_{\Phi(g)}\circ\Phi,\qquad
 \Phi\circ R_{\mathbf{g}}=R_{\Phi(g)}\circ\Phi,\qquad
 \Phi\circ C_g=C_{\Phi(g)}\circ\Phi.
\]
\end{proposition}

\begin{proof}
All identities follow directly from the group laws.
Smoothness of the maps and their inverses follows from smoothness of
multiplication and inversion.
\end{proof}

\begin{proposition}[Classical matrix groups]
\label{prop:grupos-matriciales-clasicos}
\index{general linear group}
\index{orthogonal group}
\index{unitary group}
The groups
\[
 \operatorname{GL}(n,\mathbb R),\quad
 \operatorname{GL}(n,\mathbb C),\quad
 \operatorname{O}(n),\quad
 \operatorname{SO}(n),\quad
 \operatorname{U}(n)
\]
are finite-dimensional real Lie groups with the smooth structures
induced by their inclusions into matrix spaces. In particular,
\(\operatorname{U}(1)=\{z\in\mathbb C:|z|=1\}\) is the circle viewed as a Lie
group.
\end{proposition}

\begin{proof}
The set \(\operatorname{GL}(n,\mathbb K)=\{A:\det A\neq0\}\), for
\(\mathbb K=\mathbb R,\mathbb C\), is open in the real vector space of
matrices. Multiplication is polynomial, and the formula
\(A^{-1}=\operatorname{adj}(A)/\det A\) proves that inversion is smooth.

For the orthogonal group, consider
\[
 F\colon M_n(\mathbb R)\longrightarrow\operatorname{Sym}_n(\mathbb R),
 \qquad F(A)=A^{\mathsf T}A.
\]
If \(A\in F^{-1}(I)\), then
\[
 dF_A(H)=A^{\mathsf T}H+H^{\mathsf T}A.
\]
This differential is surjective: given a symmetric matrix \(S\), simply
take \(H=AS/2\). By the regular value theorem,
\(\operatorname{O}(n)=F^{-1}(I)\) is a submanifold. The group operations
are restrictions of smooth maps. Since the determinant takes only
the values \(\pm1\) on \(\operatorname{O}(n)\), the subset
\(\operatorname{SO}(n)=\{A\in\operatorname{O}(n):\det A=1\}\) is open and
closed in \(\operatorname{O}(n)\), and is therefore also a Lie group.

The same argument, regarding complex spaces as real
spaces and using
\[
 F_{\mathbb C}(A)=A^*A,
 \qquad d(F_{\mathbb C})_A(H)=A^*H+H^*A,
\]
shows that \(\operatorname{U}(n)=F_{\mathbb C}^{-1}(I)\) is a
submanifold: for a Hermitian matrix \(S\), the choice \(H=AS/2\)
proves surjectivity. The case \(n=1\) gives precisely the unit
circle.
\end{proof}

Before defining the Lie algebra, we recall our convention for the bracket
of vector fields.

\begin{proposition}[Bracket of vector fields]
\label{prop:corchete-campos-vectoriales-naturalidad}
Let $M$ and $N$ be smooth manifolds with or without boundary, and let
$\mathbf{X},\mathbf{Y}\in\mathfrak X(M)$. There is a unique smooth vector field
$[\mathbf{X},\mathbf{Y}]$ characterized by
\begin{equation}
\label{eq:corchete-campos-como-conmutador}
 [\mathbf{X},\mathbf{Y}](f)=\mathbf{X}(\mathbf{Y}f)-\mathbf{Y}(\mathbf{X}f),
 \qquad f\in C^\infty(M).
\end{equation}
The bracket is bilinear and antisymmetric and satisfies the Jacobi identity.
If $F\colon M\to N$ is smooth, $\mathbf{X}_1,\mathbf{X}_2\in\mathfrak X(M)$ and
$\mathbf{Y}_1,\mathbf{Y}_2\in\mathfrak X(N)$, and $\mathbf{X}_i$ is $F$-related to $\mathbf{Y}_i$ for
$i=1,2$, then $[\mathbf{X}_1,\mathbf{X}_2]$ is $F$-related to $[\mathbf{Y}_1,\mathbf{Y}_2]$.
In particular, if $F$ is a diffeomorphism,
\[
 F_*[\mathbf{X},\mathbf{Y}]=[F_*\mathbf{X},F_*\mathbf{Y}].
\]
\end{proposition}

\begin{proof}
For $f,g\in C^\infty(M)$, the Leibniz rule gives
\[
 \bigl(\mathbf{X}\mathbf{Y}-\mathbf{Y}\mathbf{X}\bigr)(fg)
 =f\bigl(\mathbf{X}\mathbf{Y}-\mathbf{Y}\mathbf{X}\bigr)g
 +g\bigl(\mathbf{X}\mathbf{Y}-\mathbf{Y}\mathbf{X}\bigr)f.
\]
Thus, at each point $p$, the operator
$f\mapsto (\mathbf{X}\mathbf{Y}-\mathbf{Y}\mathbf{X})f(p)$ is a derivation and therefore a tangent vector
at $T_pM$. In coordinates, if
$\mathbf{X}=X^i\partial_i$ and $\mathbf{Y}=Y^j\partial_j$, then
\[
 [\mathbf{X},\mathbf{Y}]
 =
 \left(X^i\partial_iY^j-Y^i\partial_iX^j\right)\partial_j,
\]
so the resulting field is smooth. Antisymmetry and the Jacobi
identity are the corresponding identities for commutators of
linear operators.

If $\mathbf{X}_i$ is $F$-related to $\mathbf{Y}_i$, then, for
$f\in C^\infty(N)$,
\[
 \mathbf{X}_i(f\circ F)=(\mathbf{Y}_i f)\circ F.
\]
Applying this equality twice gives
\[
 [\mathbf{X}_1,\mathbf{X}_2](f\circ F)
 =
 \bigl([\mathbf{Y}_1,\mathbf{Y}_2]f\bigr)\circ F,
\]
which is precisely the required relation.
\end{proof}

\begin{definition}[Invariant fields and the Lie algebra]
\label{def:algebra-lie-campos-invariantes}
\index{vector field!left invariant}
\index{Lie algebra@Lie algebra!of a Lie group}
Let $G$ be a finite-dimensional real Lie group whose underlying smooth
manifold has no boundary. A field
$\mathbf{X}\in\mathfrak X(G)$ is \textit{left invariant} if
\[
 (L_g)_*\mathbf{X}=\mathbf{X}
 \qquad\text{for every }g\in G.
\]
For $\xi\in T_eG$, define
\[
 \boldsymbol{\xi}^{L}_g:=(dL_g)_e\xi.
\]
The vector space
\[
 \mathfrak g:=T_eG
\]
is called the \textit{Lie algebra of $G$}. Its bracket is defined by
\begin{equation}
\label{eq:corchete-algebra-lie-izquierdo}
 [\xi,\eta]_{\mathfrak g}:=[\boldsymbol{\xi}^{L},\boldsymbol{\eta}^{L}]_e.
\end{equation}
We omit the subscript $\mathfrak g$ when there is no risk of confusion.
\end{definition}

\begin{proposition}
\label{prop:campos-izquierdos-algebra-lie}
Let $G$ be a finite-dimensional real Lie group whose underlying smooth
manifold has no boundary. The assignment
\[
 \mathfrak g\longrightarrow
 \{\mathbf{X}\in\mathfrak X(G)\mid \mathbf{X}\text{ is left invariant}\},
 \qquad
 \xi\longmapsto\boldsymbol{\xi}^{L},
\]
is a linear isomorphism with inverse $\mathbf{X}\mapsto \mathbf{X}_e$. Moreover,
\[
 [\boldsymbol{\xi}^{L},\boldsymbol{\eta}^{L}]
 =\boldsymbol{[\xi,\eta]}^{L}.
\]
Therefore, \eqref{eq:corchete-algebra-lie-izquierdo} makes
$\mathfrak g$ a Lie algebra.
\end{proposition}

\begin{proof}
Smoothness of $\boldsymbol{\xi}^{L}$ follows from smoothness of
$(g,\xi)\mapsto(dL_g)_e\xi$. The identity
$L_g\circ L_h=L_{gh}$ implies
\[
 (dL_g)_h(\boldsymbol{\xi}^{L}_h)
 =(dL_g)_h(dL_h)_e\xi
 =(dL_{gh})_e\xi
 =\boldsymbol{\xi}^{L}_{gh},
\]
so $\boldsymbol{\xi}^{L}$ is left invariant. Conversely, if
$\mathbf{X}$ is left invariant, then
\[
 \mathbf{X}_g=(dL_g)_e\mathbf{X}_e,
\]
so $\mathbf{X}=\boldsymbol{(\mathbf{X}_e)}^{L}$. Naturality of the bracket under $L_g$ shows that the
bracket of two left invariant fields is again left
invariant. Evaluation at $e$ gives
$[\boldsymbol{\xi}^{L},\boldsymbol{\eta}^{L}]
=\boldsymbol{[\xi,\eta]}^{L}$. The Lie algebra identities are then inherited
from the bracket of vector fields.
\end{proof}

\begin{proposition}[Lie algebras of the classical matrix groups]
\label{prop:algebras-lie-grupos-matriciales-clasicos}
With the matrix bracket \([A,B]=AB-BA\), we have the identifications
\[
 \mathfrak{gl}(n,\mathbb R)=M_n(\mathbb R),\qquad
 \mathfrak{gl}(n,\mathbb C)=M_n(\mathbb C)
\]
as real Lie algebras, and
\[
 \mathfrak{o}(n)=\mathfrak{so}(n)
 =\{A\in M_n(\mathbb R):A^{\mathsf T}+A=0\},
 \qquad
 \mathfrak{u}(n)=\{A\in M_n(\mathbb C):A^*+A=0\}.
\]
In particular, \(\mathfrak u(1)=i\mathbb R\).
\end{proposition}

\begin{proof}
The tangent space at the identity of a general linear group is the entire
ambient space. If \(A(t)\) is a curve in \(\operatorname{O}(n)\) with
\(A(0)=I\), differentiating \(A(t)^{\mathsf T}A(t)=I\) gives
\(A'(0)^{\mathsf T}+A'(0)=0\). Conversely, each skew-symmetric matrix
belongs to the kernel of \(dF_I\), which is the tangent space to the regular
level set. Since \(\operatorname{SO}(n)\) is open in
\(\operatorname{O}(n)\), both groups have the same tangent space at
the identity. The identity \(A(t)^*A(t)=I\) similarly gives the description
of \(\mathfrak u(n)\).

Finally, for a matrix \(A\), the left invariant field is
\(\mathbf X_A(g)=gA\). A calculation in the matrix space gives
\([\mathbf X_A,\mathbf X_B](g)=g(AB-BA)\), so the bracket defined
in \eqref{eq:corchete-algebra-lie-izquierdo} agrees with the commutator.
\end{proof}

\begin{proposition}[Differential of a homomorphism]
\label{prop:diferencial-homomorfismo-lie}
Let $G$ and $H$ be finite-dimensional real Lie groups whose underlying smooth
manifolds have no boundary, and let
$\Phi\colon G\to H$ be a Lie group homomorphism. Then
\[
 d\Phi_e\colon\mathfrak g\longrightarrow\mathfrak h
\]
is a Lie algebra homomorphism:
\[
 d\Phi_e[\xi,\eta]=[d\Phi_e\xi,d\Phi_e\eta].
\]
\end{proposition}

\begin{proof}
The identity
$\Phi\circ L_g=L_{\Phi(g)}\circ\Phi$ implies
\[
 d\Phi_g(\boldsymbol{\xi}^{L}_g)
 =
 (dL_{\Phi(g)})_e(d\Phi_e\xi)
 =
 \boldsymbol{(d\Phi_e\xi)}^L_{\Phi(g)}.
\]
Thus $\boldsymbol{\xi}^{L}$ and $\boldsymbol{(d\Phi_e\xi)}^L$ are
$\Phi$-related.
Proposition~\ref{prop:corchete-campos-vectoriales-naturalidad} shows that
their brackets are also related. Evaluating at $e$ gives the formula.
\end{proof}

\begin{lemma}[Bracket and flow]
\label{lem:corchete-mediante-flujo}
Let $M$ be a smooth manifold without boundary, let
$\mathbf{X}\in\mathfrak X(M)$, and let $\Phi_t$ be its local flow. For every
$\mathbf{Y}\in\mathfrak X(M)$,
\begin{equation}
\label{eq:corchete-mediante-flujo}
 [\mathbf{X},\mathbf{Y}]
 =
 \left.\frac{d}{dt}\right|_{t=0}(\Phi_{-t})_*\mathbf{Y},
\end{equation}
on any open set where both sides are defined.
\end{lemma}

\begin{proof}
As operators on functions,
\[
 \bigl((\Phi_{-t})_*\mathbf{Y}\bigr)f
 =
 \Phi_t^*\bigl(\mathbf{Y}(\Phi_{-t}^*f)\bigr).
\]
The derivative at $t=0$ of the operator $\Phi_t^*$ is $\mathbf{X}$, while the
derivative of $\Phi_{-t}^*$ is $-\mathbf{X}$. By the product rule,
\[
 \left.\frac{d}{dt}\right|_{0}
 \Phi_t^*\mathbf{Y}\Phi_{-t}^*
 =
 \mathbf{X}\mathbf{Y}-\mathbf{Y}\mathbf{X}=[\mathbf{X},\mathbf{Y}].
\]
\end{proof}

\begin{theorem}[Exponential map of a Lie group]
\label{teo:exponencial-grupo-lie}
\index{exponential map!of a Lie group}
Let $G$ be a finite-dimensional real Lie group whose underlying smooth
manifold has no boundary. For each $\xi\in\mathfrak g$, the field
$\boldsymbol{\xi}^{L}$ is complete. If
$\gamma_\xi\colon\mathbb R\to G$ is its integral curve with
$\gamma_\xi(0)=e$, then
\[
 \gamma_\xi(t+s)=\gamma_\xi(t)\gamma_\xi(s).
\]
The map
\[
 \exp_G\colon\mathfrak g\longrightarrow G,\qquad
 \exp_G(\xi):=\gamma_\xi(1),
\]
is smooth and satisfies
\begin{equation}
\label{eq:propiedades-basicas-exponencial-lie}
 \exp_G(t\xi)=\gamma_\xi(t),\qquad
 \exp_G((t+s)\xi)=\exp_G(t\xi)\exp_G(s\xi).
\end{equation}
Moreover,
\[
 (d\exp_G)_0=\operatorname{Id}_{\mathfrak g};
\]
consequently, $\exp_G$ restricts to a diffeomorphism between an open
neighborhood of $0\in\mathfrak g$ and an open neighborhood of $e\in G$.
\end{theorem}

\begin{proof}
Let $\gamma_\xi\colon I\to G$ be the maximal integral curve of $\boldsymbol{\xi}^{L}$
starting at $e$. If $s\in I$, the curves
\[
 t\longmapsto\gamma_\xi(s+t),
 \qquad
 t\longmapsto\gamma_\xi(s)\gamma_\xi(t)
\]
are integral curves of $\boldsymbol{\xi}^{L}$ and take the same value at $t=0$. By
uniqueness, they agree wherever both are defined.

The interval $I$ contains some $(-\varepsilon,\varepsilon)$. If
$b:=\displaystyle\sup I<\infty$, choose $s\in I$ with $b-s<\varepsilon$. The formula
\[
 \widetilde\gamma(s+t):=\gamma_\xi(s)\gamma_\xi(t),
 \qquad |t|<\varepsilon,
\]
agrees with $\gamma_\xi$ on the overlap and extends the curve beyond
$b$, contradicting maximality. The same argument at the left
endpoint shows that $I=\mathbb R$. The group identity already obtained
therefore holds for every $s,t\in\mathbb R$.

The rescaling lemma for integral curves gives
\[
 \gamma_{t\xi}(1)=\gamma_\xi(t),
\]
which proves \eqref{eq:propiedades-basicas-exponencial-lie}. To prove
smoothness with respect to $\xi$, consider the smooth field on
$G\times\mathfrak g$
\[
 \boldsymbol{\mathcal X}_{(g,\xi)}:=\bigl((dL_g)_e\xi,0\bigr).
\]
Each integral curve of $\boldsymbol{\mathcal X}$ keeps its second component fixed, and
its first component is an integral curve of a left invariant
field; by what we have proved, $\boldsymbol{\mathcal X}$ is complete. The fundamental theorem
of flows implies that its global flow depends smoothly on the
initial point. Evaluating it at time $1$ and at $(e,\xi)$ shows that
$\xi\mapsto\exp_G(\xi)$ is smooth.

Finally,
\[
 (d\exp_G)_0(\xi)
 =
 \left.\frac{d}{dt}\right|_{0}\exp_G(t\xi)
 =
 \gamma_\xi'(0)=\xi.
\]
Theorem~\ref{teo: funcion inversa variedades} gives the local
assertion.
\end{proof}

\begin{corollary}[One-parameter subgroups]
\label{cor:subgrupos-un-parametro-lie}
\index{one-parameter subgroup@one-parameter subgroup}
Let $G$ be a finite-dimensional real Lie group whose underlying smooth
manifold has no boundary. For each $\xi\in\mathfrak g$, the map
\[
 t\longmapsto\exp_G(t\xi)
\]
is a smooth homomorphism $(\mathbb R,+)\to G$. Conversely, if
$\alpha\colon(\mathbb R,+)\to G$ is a smooth homomorphism and
$\xi=\alpha'(0)$, then
\[
 \alpha(t)=\exp_G(t\xi)
 \qquad\text{for every }t\in\mathbb R.
\]
\end{corollary}

\begin{proof}
The first assertion is \eqref{eq:propiedades-basicas-exponencial-lie}.
For the converse,
\[
 \alpha(t+s)=\alpha(t)\alpha(s)
\]
implies, on differentiating with respect to $s$ at $0$,
\[
 \alpha'(t)=(dL_{\alpha(t)})_e\xi=\boldsymbol{\xi}^{L}_{\alpha(t)}.
\]
Thus $\alpha$ and $t\mapsto\exp_G(t\xi)$ are integral curves of
$\boldsymbol{\xi}^{L}$ with
the same initial data, and uniqueness proves the equality.
\end{proof}

\begin{definition}[Adjoint and infinitesimal representations]
\label{def:representaciones-ad-ad-lie}
\index{adjoint representation@adjoint representation}
Let $G$ be a finite-dimensional real Lie group whose underlying smooth
manifold has no boundary. Define
\[
 \operatorname{Ad}_g:=(dC_g)_e\in GL(\mathfrak g),
 \qquad
 \operatorname{Ad}\colon G\longrightarrow GL(\mathfrak g),
\]
and
\[
 \operatorname{ad}_\xi\eta
 :=
 \left.\frac{d}{dt}\right|_{0}
 \operatorname{Ad}_{\exp_G(t\xi)}\eta.
\]
\end{definition}

\begin{proposition}[Naturality and equivariance of the exponential map]
\label{prop:naturalidad-equivariancia-exponencial-lie}
\index{exponential map!naturality}\index{exponential map!equivariance}
Let $G$ and $H$ be finite-dimensional real Lie groups whose underlying smooth
manifolds have no boundary. If
$\Phi\colon G\to H$ is a Lie group homomorphism, then
\[
 \Phi(\exp_G\xi)=\exp_H(d\Phi_e\xi),
 \qquad \xi\in\mathfrak g.
\]
In particular, for $g\in G$,
\begin{equation}
 \exp_G(\operatorname{Ad}_g\xi)
 =g\exp_G(\xi)g^{-1}.
 \label{eq:equivariancia-exponencial-ad}
\end{equation}
\end{proposition}
\begin{proof}
The two maps
\[
 t\longmapsto\Phi(\exp_G(t\xi)),
 \qquad
 t\longmapsto\exp_H(t\,d\Phi_e\xi)
\]
are one-parameter subgroups of $H$ with initial velocity
$d\Phi_e\xi$. Uniqueness in
Corollary~\ref{cor:subgrupos-un-parametro-lie} shows that they agree;
evaluate at $t=1$. For $\Phi=C_g$, its differential at $e$ is
$\operatorname{Ad}_g$, giving
\eqref{eq:equivariancia-exponencial-ad}.
\end{proof}

\begin{proposition}[Identities for \texorpdfstring{$\operatorname{Ad}$ and
$\operatorname{ad}$}{Ad and ad}]
\label{prop:identidades-ad-ad-lie}
Let $G$ be a finite-dimensional real Lie group whose underlying smooth
manifold has no boundary. The following identities hold:
\begin{align}
 \operatorname{Ad}_{gh}
 &=\operatorname{Ad}_g\operatorname{Ad}_h,
 \label{eq:Ad-representacion}\\
 \operatorname{Ad}_g[\xi,\eta]
 &=[\operatorname{Ad}_g\xi,\operatorname{Ad}_g\eta],
 \label{eq:Ad-preserva-corchete}\\
 \operatorname{ad}_\xi\eta
 &=[\xi,\eta],
 \label{eq:ad-es-corchete}\\
 \operatorname{Ad}_{\exp_G(t\xi)}
 &=\exp\bigl(t\operatorname{ad}_\xi\bigr),
 \label{eq:Ad-exp-ad}\\
 \operatorname{Ad}_g\operatorname{ad}_\xi\operatorname{Ad}_{g^{-1}}
 &=\operatorname{ad}_{\operatorname{Ad}_g\xi},
 \label{eq:Ad-conjuga-ad}\\
 \operatorname{ad}_{[\xi,\eta]}
 &=[\operatorname{ad}_\xi,\operatorname{ad}_\eta].
 \label{eq:ad-homomorfismo}
\end{align}
\end{proposition}

\begin{proof}
The identity $C_{gh}=C_g\circ C_h$ implies
\eqref{eq:Ad-representacion}. Since $C_g$ is a Lie group automorphism,
Proposition~\ref{prop:diferencial-homomorfismo-lie} gives
\eqref{eq:Ad-preserva-corchete}.

The flow of $\boldsymbol{\xi}^{L}$ is
\[
 \Phi_t=R_{\exp_G(t\xi)},
\]
since the curve starting at $g$ is $g\exp_G(t\xi)$. A direct verification
gives
\[
 (R_h)_*\boldsymbol{\eta}^{L}
 =\boldsymbol{(\operatorname{Ad}_{h^{-1}}\eta)}^L.
\]
By Lemma~\ref{lem:corchete-mediante-flujo},
\[
 [\boldsymbol{\xi}^{L},\boldsymbol{\eta}^{L}]
 =
 \left.\frac{d}{dt}\right|_0
 (R_{\exp_G(-t\xi)})_*\boldsymbol{\eta}^{L}
 =
 \boldsymbol{\left(
 \left.\frac{d}{dt}\right|_0
 \operatorname{Ad}_{\exp_G(t\xi)}\eta
 \right)}^L.
\]
Evaluating at $e$ gives \eqref{eq:ad-es-corchete}.

By \eqref{eq:Ad-representacion}, the map
$B(t):=\operatorname{Ad}_{\exp_G(t\xi)}$ satisfies
\[
 B(t+s)=B(t)B(s),\qquad B(0)=I.
\]
Differentiating with respect to $s$ at $0$ gives
$B'(t)=B(t)\operatorname{ad}_\xi$, whose unique solution is
\eqref{eq:Ad-exp-ad}. Identity \eqref{eq:Ad-conjuga-ad} follows
by differentiating
\[
 \operatorname{Ad}_g\operatorname{Ad}_{\exp(t\xi)}
 \operatorname{Ad}_{g^{-1}}
 =
 \operatorname{Ad}_{\exp(t\operatorname{Ad}_g\xi)}.
\]
Finally, \eqref{eq:ad-homomorfismo} is the Jacobi identity written
as an identity of endomorphisms.
\end{proof}

\begin{theorem}[Closed subgroup theorem]
\label{teo:subgrupo-cerrado-lie}
\index{closed subgroup theorem}\index{Lie subgroup!closed}
Let $G$ be a finite-dimensional real Lie group whose underlying smooth
manifold has no boundary, and let $H\subseteq G$ be a subgroup closed
in the topology of $G$. Then $H$ has a unique smooth structure
making it an embedded Lie subgroup of $G$. Its Lie algebra,
identified with a subspace of $\mathfrak g$, is
\begin{equation}
 \mathfrak h
 =
 \{\xi\in\mathfrak g\mid
   \exp_G(t\xi)\in H\text{ for every }t\in\mathbb R\}.
 \label{eq:algebra-lie-subgrupo-cerrado-exponencial}
\end{equation}
\end{theorem}

Existence and the characterization of the Lie algebra are
\cite[Theorem~20.12 and Proposition~20.9, pp.~521--523]{LeeS}.
\begin{proof}
The existence assertion is the cited closed subgroup theorem. The structure is unique because an
embedded Lie subgroup structure must be the submanifold
structure induced by slice charts. Once it is constructed,
\eqref{eq:algebra-lie-subgrupo-cerrado-exponencial} is the characterization
given by the cited exponential characterization: the curves
$t\mapsto\exp_G(t\xi)$ contained in $H$ are precisely the one-parameter
subgroups whose initial velocity belongs to $T_eH$.
\end{proof}

\begin{proposition}[Centralizers]
\label{prop:centralizador-subgrupo-cerrado-lie}
\index{centralizer}\index{Lie algebra@Lie algebra!of a centralizer}
Let $G$ be a finite-dimensional real Lie group whose underlying smooth
manifold has no boundary. For any subset $S\subseteq G$,
its centralizer
\[
 Z_G(S):=\{g\in G\mid gs=sg\text{ for every }s\in S\}
\]
is a closed subgroup, and hence an embedded Lie subgroup. If $G$
is compact, $Z_G(S)$ is compact. Its Lie algebra is
\begin{equation}
 \operatorname{Lie}(Z_G(S))
 =
 \{\xi\in\mathfrak g\mid
   \operatorname{Ad}_s\xi=\xi\text{ for every }s\in S\}.
 \label{eq:algebra-lie-centralizador-fijos-ad}
\end{equation}
If $S=H$ is a connected Lie subgroup with algebra $\mathfrak h$, this
expression is equivalent to
\[
 \operatorname{Lie}(Z_G(H))
 =
 \{\xi\in\mathfrak g\mid[\xi,\eta]=0
   \text{ for every }\eta\in\mathfrak h\}.
\]
\end{proposition}
\begin{proof}
The centralizer is a subgroup by the group laws, and
\[
 Z_G(S)=\bigcap_{s\in S}
 \{g\in G\mid C_g(s)=s\}.
\]
Each set on the right is closed, since $g\mapsto C_g(s)$ is
continuous and $\{s\}$ is closed.
Theorem~\ref{teo:subgrupo-cerrado-lie} gives it its embedded structure; if
$G$ is compact, the closed subset $Z_G(S)$ is also compact.

By \eqref{eq:algebra-lie-subgrupo-cerrado-exponencial}, a vector
$\xi$ belongs to its algebra if and only if
$\exp_G(t\xi)$ commutes with every $s\in S$ for every $t$. Using
\eqref{eq:equivariancia-exponencial-ad}, this condition becomes
\[
 \exp_G(t\operatorname{Ad}_s\xi)=\exp_G(t\xi)
 \quad(s\in S,\ t\in\mathbb R).
\]
Differentiating at $t=0$ gives $\operatorname{Ad}_s\xi=\xi$; the converse follows
from the same equivariance formula. This proves
\eqref{eq:algebra-lie-centralizador-fijos-ad}.

If $H$ is connected, it is generated by any neighborhood of the identity, and such a neighborhood
can be covered by exponentials. By
$\operatorname{Ad}_{\exp(t\eta)}
=\exp(t\operatorname{ad}_\eta)$, a vector is fixed by every $H$ if and
only if $\operatorname{ad}_\eta\xi=[\eta,\xi]=0$ for every
$\eta\in\mathfrak h$.
\end{proof}

\begin{definition}[Smooth actions and positive fundamental fields]
\label{def:acciones-suaves-campos-fundamentales-positivos}
\index{smooth action@smooth action}
\index{fundamental field!positive}
Let $G$ be a finite-dimensional real Lie group whose underlying smooth
manifold has no boundary, and let $M$ be a smooth manifold with or without
boundary. A \textit{smooth left action} of $G$ on $M$ is
a smooth map
\[
 a\colon G\times M\longrightarrow M,\qquad (g,x)\longmapsto g\cdot x,
\]
such that
\[
 e\cdot x=x,\qquad (gh)\cdot x=g\cdot(h\cdot x).
\]
A \textit{smooth right action} is a smooth map
\[
 r\colon M\times G\longrightarrow M,\qquad (x,g)\longmapsto x\cdot g,
\]
such that
\[
 x\cdot e=x,\qquad (x\cdot g)\cdot h=x\cdot(gh).
\]

We always use the \textit{positive fundamental field}. For a left
action, define
\[
 \boldsymbol{\xi}_M(x)
 :=
 \left.\frac{d}{dt}\right|_0\exp_G(t\xi)\cdot x,
\]
and for a right action, define
\[
 \boldsymbol{\xi}_M(x)
 :=
 \left.\frac{d}{dt}\right|_0x\cdot\exp_G(t\xi).
\]
\end{definition}

\begin{proposition}[Signs and transformation of fundamental fields]
\label{prop:signos-transformacion-campos-fundamentales}
Let $G$ be a finite-dimensional real Lie group whose underlying smooth
manifold has no boundary, let $M$ be a smooth manifold with or without
boundary, consider a smooth left or right action of $G$ on
$M$, and let $\xi,\eta\in\mathfrak g$.

\begin{enumerate}[label=(\alph*)]
\item For a left action, the flow of $\boldsymbol{\xi}_M$ is
\[
 \Phi_t(x)=\exp_G(t\xi)\cdot x,
\]
and
\begin{align}
 [\boldsymbol{\xi}_M,\boldsymbol{\eta}_M]
 &=-\boldsymbol{[\xi,\eta]}_M,
 \label{eq:signo-campos-fundamentales-izquierda}\\
 (a_g)_*\boldsymbol{\xi}_M
 &=\boldsymbol{(\operatorname{Ad}_g\xi)}_M,
 \qquad a_g(x):=g\cdot x.
 \label{eq:transformacion-fundamental-izquierda}
\end{align}

\item For a right action, the flow of $\boldsymbol{\xi}_M$ is
\[
 \Psi_t(x)=x\cdot\exp_G(t\xi),
\]
and
\begin{align}
 [\boldsymbol{\xi}_M,\boldsymbol{\eta}_M]
 &=\boldsymbol{[\xi,\eta]}_M,
 \label{eq:signo-campos-fundamentales-derecha}\\
 (r_g)_*\boldsymbol{\xi}_M
 &=\boldsymbol{(\operatorname{Ad}_{g^{-1}}\xi)}_M,
 \qquad r_g(x):=x\cdot g.
 \label{eq:transformacion-fundamental-derecha}
\end{align}
\end{enumerate}
\end{proposition}

\begin{proof}
The flow formulas follow immediately from the action law and
the group law for $t\mapsto\exp_G(t\xi)$.

For a left action,
\[
 g\exp(t\xi)\cdot x
 =
 \exp(t\operatorname{Ad}_g\xi)\cdot(g\cdot x),
\]
and differentiating at $t=0$ gives
\eqref{eq:transformacion-fundamental-izquierda}. Applying this identity with
$g=\exp(-t\xi)$ and using
Lemma~\ref{lem:corchete-mediante-flujo},
\[
 [\boldsymbol{\xi}_M,\boldsymbol{\eta}_M]
 =
 \left.\frac{d}{dt}\right|_0
 (a_{\exp(-t\xi)})_*\boldsymbol{\eta}_M
 =
 \left.\frac{d}{dt}\right|_0
 \boldsymbol{(\operatorname{Ad}_{\exp(-t\xi)}\eta)}_M
 =
 -\boldsymbol{[\xi,\eta]}_M.
\]

For a right action,
\[
 (x\cdot\exp(t\xi))\cdot g
 =
 (x\cdot g)\cdot
 \exp(t\operatorname{Ad}_{g^{-1}}\xi),
\]
which proves \eqref{eq:transformacion-fundamental-derecha}. Applying it to
$g=\exp(-t\xi)$,
\[
 [\boldsymbol{\xi}_M,\boldsymbol{\eta}_M]
 =
 \left.\frac{d}{dt}\right|_0
 (r_{\exp(-t\xi)})_*\boldsymbol{\eta}_M
 =
 \left.\frac{d}{dt}\right|_0
 \boldsymbol{(\operatorname{Ad}_{\exp(t\xi)}\eta)}_M
 =
 \boldsymbol{[\xi,\eta]}_M.
\]
\end{proof}

\begin{definition}[Orbit, stabilizer, and orbit map]
\label{def:orbita-estabilizador-accion-lie}
Let $G$ be a finite-dimensional real Lie group whose underlying smooth
manifold has no boundary, let $M$ be a smooth manifold with or without
boundary, and consider a smooth left or right action of $G$ on
$M$. For a left action and $x\in M$, define
\[
 \mathcal O_x:=G\cdot x,
 \qquad
 G_x:=\{g\in G\mid g\cdot x=x\},
 \qquad
 \vartheta_x\colon G\to M,\qquad \vartheta_x(g)=g\cdot x.
\]
For a right action, use
$\mathcal O_x=x\cdot G$, $G_x=\{g\mid x\cdot g=x\}$, and
$\vartheta_x(g)=x\cdot g$. An action is \emph{free} if all its
stabilizers are trivial, and \emph{transitive} if it has a single
orbit.
\end{definition}

\begin{proposition}[Elementary properties of orbits]
\label{prop:orbitas-estabilizadores-campos-fundamentales}
Let $G$ be a finite-dimensional real Lie group whose underlying smooth
manifold has no boundary, let $M$ be a smooth manifold with or without
boundary, and consider a smooth left or right action of $G$ on
$M$. The orbit map is smooth, $G_x$ is a closed embedded Lie subgroup
of $G$, and
\begin{equation}
\label{eq:diferencial-mapeo-orbita-campo-fundamental}
 (d\vartheta_x)_e\xi=\boldsymbol{\xi}_M(x).
\end{equation}
For a left action,
$G_{g\cdot x}=gG_xg^{-1}$; for a right action,
$G_{x\cdot g}=g^{-1}G_xg$. In particular, freeness and the dimension of the
kernel of $(d\vartheta_x)_e$ are constant along each orbit.
\end{proposition}

\begin{proof}
Smoothness follows from that of the action with $x$ fixed. The action laws show
that $G_x$ contains $e$ and is closed under products and inverses. Moreover,
$G_x=\vartheta_x^{-1}(\{x\})$ is closed because manifolds are
Hausdorff; Theorem~\ref{teo:subgrupo-cerrado-lie} gives it its embedded Lie subgroup
structure. The identity
\eqref{eq:diferencial-mapeo-orbita-campo-fundamental} is exactly the
definition of the fundamental field. Finally,
\[
 h\cdot(g\cdot x)=g\cdot x
 \iff g^{-1}hg\cdot x=x
\]
proves the formula for a left action; the right-action formula follows from
$(x\cdot g)\cdot h=x\cdot g$ by multiplying on the right by $g^{-1}$.
The transformation of the kernels also follows from
Proposition~\ref{prop:signos-transformacion-campos-fundamentales}.
\end{proof}

\begin{theorem}[Haar measure on compact groups]
\label{teo:medida-haar-compacto}
\index{Haar measure!compact group}
Let $G$ be a compact finite-dimensional real Lie group whose underlying smooth
manifold has no boundary. There exists a unique regular probability measure
$\mu_G$ on $G$ invariant under left and right
translations:
\[
 (L_g)_*\mu_G=\mu_G=(R_g)_*\mu_G
 \qquad(g\in G).
\]
For the construction of the normalized Haar volume form on a compact Lie
group, see \cite[Proposition~16.10]{LeeS}. Existence and
uniqueness of Haar measure is the external result used here; right
invariance follows by applying uniqueness of the left invariant probability
measure to each right translation.
\end{theorem}

\begin{proposition}[Averaging an inner product]
\label{prop:promedio-representacion-compacta}
\index{averaging over a compact group}
Let $G$ be a compact finite-dimensional real Lie group whose underlying smooth
manifold has no boundary, and let
$\rho\colon G\to GL(V)$ be a smooth representation on a finite-dimensional real
vector space. Then $V$ admits a
$G$-invariant inner product.
\end{proposition}
\begin{proof}
Start with any inner product $\langle\cdot,\cdot\rangle_0$ and
define
\[
 \langle v,w\rangle_G
 :=
 \int_G\langle\rho(g)v,\rho(g)w\rangle_0\,d\mu_G(g).
\]
The resulting form is bilinear, symmetric, and positive definite. For
$h\in G$, right invariance of $\mu_G$ and
$\rho(gh)=\rho(g)\rho(h)$ give
\[
 \langle\rho(h)v,\rho(h)w\rangle_G
 =
 \int_G\langle\rho(gh)v,\rho(gh)w\rangle_0\,d\mu_G(g)
 =\langle v,w\rangle_G.
\]
\end{proof}

\begin{corollary}[\texorpdfstring{$\operatorname{Ad}$}{Ad}-invariant
inner product]
\label{cor:producto-interno-Ad-invariante-compacto}
\index{Ad-invariant inner product@\(\operatorname{Ad}\)-invariant
inner product}
Let $G$ be a compact finite-dimensional real Lie group whose underlying smooth
manifold has no boundary. Its Lie algebra $\mathfrak g$ admits an inner
product $\langle\cdot,\cdot\rangle_{\mathfrak g}$ such that
\[
 \langle\operatorname{Ad}_g\xi,\operatorname{Ad}_g\eta\rangle_{\mathfrak g}
 =
 \langle\xi,\eta\rangle_{\mathfrak g}
 \qquad(g\in G).
\]
\end{corollary}

\begin{proof}
Apply Proposition~\ref{prop:promedio-representacion-compacta} to the
representation $\rho=\operatorname{Ad}$.
\end{proof}

\section{Tensor fields and differential forms}

Tangent and cotangent spaces allow us to construct tensors pointwise. Assembling them over $M$ gives tensor bundles, whose sections describe metrics, differential forms, curvature, and many other geometric fields. This final section of the chapter first recalls the necessary multilinear algebra and then organizes it smoothly over the manifold.

 \begin{definition}\label{def:nociones-fundamentales-multilineal}\index{multilinear}
 If $V_{1},\dots, V_{k}$ and $W$ are vector spaces, a function $F\colon V_{1}\times\dots\times V_{k}\longrightarrow W$ is called \textbf{multilinear} if it is linear in each argument when the others are fixed. The space of multilinear functions is denoted by $L(V_{1},\dots,V_{k};W)$.
 \end{definition}
 \begin{definition}\label{def:nociones-fundamentales-producto-tensorial-de-con}\index{tensor product}
 If $V_{1},\dots, V_{k}$ and $W_{1},\dots, W_{l}$ are vector spaces, and $\omega\in L(V_{1},\dots,V_{k};\mathbb{R})$ and $\eta\in L(W_{1},\dots,W_{l};\mathbb{R})$, define a function $\omega\otimes\eta \colon V_{1}\times\dots\times V_{k}\times W_{1}\times\dots\times W_{l}\longrightarrow \mathbb{R}$, called the \textbf{tensor product of $\omega$ and $\eta$}, by \[\omega\otimes\eta(v_{1},\dots,v_{k},w_{1},\dots,w_{l})=\omega(v_{1},\dots,v_{k})\eta(w_{1},\dots,w_{l}).\]
 \end{definition}
 The space of multilinear functions is a vector space, and if the underlying vector spaces are finite dimensional, it is also finite dimensional:
 \begin{proposition}\label{prop:nociones-fundamentales-espacios-vectoriales-reales-espacio-vectorial-real}
 Let $V_{1},\dots,V_{k}$ and $W$ be real vector spaces. Then $L(V_{1},\dots,V_{k};W)$ is a real vector space under the operations \[(F+G)(v_{1},\dots,v_{k})=F(v_{1},\dots,v_{k})+G(v_{1},\dots,v_{k}),\] \[(aF)(v_{1},\dots,v_{k})=a(F(v_{1},\dots,v_{k}))\] for each $F,G\in L(V_{1},\dots,V_{k};W)$ and each $a\in\mathbb{R}$.
 \end{proposition}
\begin{proof}
Fixing all variables except the $j$th, the functions
$F+G$ and $aF$ are, respectively, a sum and a scalar multiple of
linear maps $V_j\to W$; hence they are again linear in that
variable. Since this holds for each $j$, both functions are multilinear.
The vector space identities are checked pointwise and reduce
to the corresponding identities in $W$. The zero function is the identity
element, and the additive inverse of $F$ is the function $-F$.
\end{proof}

 \begin{proposition}\label{prop:nociones-fundamentales-espacios-vectoriales-dimensiones-respectivamente-base-corres}
 Let $V_{1},\dots,V_{k}$ be vector spaces of dimensions $n_{1},\dots,n_{k}$, respectively. For each $j\in \{1,\dots,k\},$, let $\{\mathbf{E}_{1}^{(j)},\dots,\mathbf{E}_{n_{j}}^{(j)}\}$ be a basis for $V_{j}$ and let $\{\varepsilon_{(j)}^{1},\dots,\varepsilon_{(j)}^{n_{j}}\}$ be the corresponding dual basis for $V_{j}^{*}$. Then \[\mathscr{B}=\{\varepsilon_{(1)}^{i_{1}}\otimes\dots\otimes\varepsilon_{(k)}^{i_{k}}\mid 1\leq i_{1}\leq n_{1},\dots,1\leq i_{k}\leq n_{k}\}\] is a basis for $L(V_{1},\dots,V_{k};\mathbb{R})$, so this vector space has dimension $n_{1}\cdots n_{k}$.
 \end{proposition}
\begin{proof}
Let $F\in L(V_1,\dots,V_k;\mathbb R)$ and write
\[
v_j=\displaystyle\sum_{i_j=1}^{n_j}v_j^{i_j}\mathbf E_{i_j}^{(j)}
\qquad (1\leq j\leq k).
\]
Successive application of linearity in each argument gives
\[
F(v_1,\dots,v_k)
=
\displaystyle\sum_{i_1=1}^{n_1}\cdots\displaystyle\sum_{i_k=1}^{n_k}
F(\mathbf E_{i_1}^{(1)},\dots,\mathbf E_{i_k}^{(k)})
v_1^{i_1}\cdots v_k^{i_k}.
\]
Since $v_j^{i_j}=\varepsilon_{(j)}^{i_j}(v_j)$, this equality expresses
$F$ as a linear combination of the elements of $\mathscr B$. To prove
independence, suppose a combination with coefficients
$c_{i_1\dots i_k}$ is zero. Evaluating it at
$(\mathbf E_{r_1}^{(1)},\dots,\mathbf E_{r_k}^{(k)})$, the identities
$\varepsilon_{(j)}^{i_j}(\mathbf E_{r_j}^{(j)})=\delta^{i_j}_{r_j}$
give $c_{r_1\dots r_k}=0$. All coefficients vanish, so
$\mathscr B$ is a basis, and its cardinality is
$n_1\cdots n_k$.
\end{proof}
 Finally, the preceding definition of the tensor product of two multilinear functions allows us to define the tensor product of vector spaces as a vector space satisfying the following \emph{universal property}:
 \begin{definition}[Universal property of the tensor product]\label{def:nociones-fundamentales-propiedad-universal-del-producto-tensorial}\index{universal property of the tensor product}
 Let $V_{1},\dots,V_{k}$ be vector spaces.
 The \textbf{tensor product} $V_{1}\otimes\dots\otimes V_{k}$ is a
 vector space equipped with a canonical multilinear map
 $\pi\colon V_1\times\cdots\times V_k\to
 V_1\otimes\cdots\otimes V_k$ having the following property: if
 $F\colon V_{1}\times\dots\times V_{k}\longrightarrow X$ is any
 multilinear function taking values in a vector space $X$, there
 exists a unique linear function
 $\widetilde{F}\colon V_{1}\otimes\dots\otimes V_{k}\longrightarrow X$
 such that the following diagram commutes:
 \[
\begin{tikzcd}
 V_{1}\times \dots\times V_{k} \arrow{r}{F}\arrow{d}{\pi} & X \\
 V_{1}\otimes\dots\otimes V_{k} \arrow[dashed,swap]{ur}{\exists !\widetilde{F}} & {}
\end{tikzcd}
\]
 \end{definition}

\begin{remark}[Existence and uniqueness of the tensor product]
Let $\mathbb R^{(V_1\times\cdots\times V_k)}$ be the free vector space
generated by the symbols $[v_1,\dots,v_k]$. Take the quotient by
the subspace generated, for each argument, by the relations
\[
[v_1,\dots,av_j+bw_j,\dots,v_k]
-a[v_1,\dots,v_j,\dots,v_k]
-b[v_1,\dots,w_j,\dots,v_k].
\]
The class of $[v_1,\dots,v_k]$ defines the map $\pi$, and the quotient
relations ensure that it is multilinear. If $F$ is multilinear, the
assignment $[v_1,\dots,v_k]\mapsto F(v_1,\dots,v_k)$ vanishes on the
subspace of relations and therefore induces a unique linear map
$\widetilde F$. This proves existence. If two pairs
$(T,\pi_T)$ and $(S,\pi_S)$ satisfy the universal property, it yields
unique linear maps $T\to S$ and $S\to T$ preserving simple
tensors; their compositions are the identity by the same uniqueness. Thus the
tensor product is unique up to a unique isomorphism compatible with the
canonical maps.
\end{remark}

 The following proposition relates the space of multilinear functions to the tensor product of vector spaces:
 \begin{proposition}\label{prop:nociones-fundamentales-espacios-vectoriales-hay-isomorfismo-canonico}
 Let $V_{1},\dots,V_{k}$ be finite-dimensional real vector spaces.
 Then there is a canonical isomorphism
 \[
 V_{1}^{*}\otimes\dots\otimes V_{k}^{*}
 \cong L(V_{1},\dots,V_{k};\mathbb{R}).
 \]
 \end{proposition}
\begin{proof}
The map
\[
(\lambda_1,\dots,\lambda_k)
\longmapsto
\bigl[(v_1,\dots,v_k)\longmapsto
\lambda_1(v_1)\cdots\lambda_k(v_k)\bigr]
\]
is multilinear in $\lambda_1,\dots,\lambda_k$. The universal property
therefore induces a canonical linear map
\[
\Phi\colon V_1^*\otimes\cdots\otimes V_k^*
\longrightarrow L(V_1,\dots,V_k;\mathbb R).
\]
In the dual bases of the preceding proposition, $\Phi$ sends the tensor
$\varepsilon_{(1)}^{i_1}\otimes\cdots\otimes
\varepsilon_{(k)}^{i_k}$ to the element of the basis $\mathscr B$ with the
same indices. Thus $\Phi$ sends a basis to a basis and is an
isomorphism. The construction is independent of the bases, since $\Phi$ was
defined solely by evaluation.
\end{proof}
 We define tensors on a given vector space as follows:
\begin{definition}\label{def:nociones-fundamentales-espacio-vectorial-dimension-finita-espacio-vectorial}\index{tensor!mixed on a vector space}
 Let $V$ be a finite-dimensional vector space. Define the vector space of mixed tensors of type $(k,l)$, $k$ times contravariant and $l$ times covariant, on $V$ by
 \[
 T^{(k,l)}(V):=\underbrace{V\otimes \dots \otimes V}_{\text{$k$ times}}\otimes \underbrace{V^{*}\otimes \dots \otimes V^{*}}_{\text{$l$ times}}.
 \]
\end{definition}

We are ready to define tensor bundles on a smooth manifold:

\begin{definition}\label{def:nociones-fundamentales-haz-de-tensores-mixtos-de-tipo}\index{tensor bundle!mixed}
 Let $M$ be a smooth manifold with or without boundary. Define the \textbf{bundle of mixed tensors of type $(k,l)$} by
 \[
 T^{(k,l)}(TM):=\coprod_{p\in M}T^{(k,l)}(T_{p}M).
 \]
\end{definition}

\begin{remark}\label{obs:nociones-fundamentales-podemos-hacer-siguientes-identificaciones-label}
 We have the following identifications:
 \begin{enumerate}[label=(\alph*)]
 \item $T^{(0,0)}(TM)=M\times \mathbb{R}$.
 \item $T^{(0,1)}(TM)=T^{*}M$.
 \item $T^{(1,0)}(TM)=TM$.
 \end{enumerate}
\end{remark}

 We can equip tensor bundles with a smooth structure as follows:
 \begin{proposition}\label{estructura del haz tensorial}
 Let $M$ be a smooth manifold with or without boundary and let $\dim M = n$.
 Then, for every pair of integers $k,l \geq 0$, the tensor bundle
 \[
 T^{(k,l)}(TM)\longrightarrow M
 \]
 admits a unique topology and a unique smooth structure making it a vector bundle of rank $n^{k+l}$ over $M$ for which the component trivializations induced by smooth charts of $M$ are smooth.
\end{proposition}

 \begin{proof}
 We apply Lemma~\ref{lema del haz vectorial suave}, first constructing the maps that will serve as local trivializations.

 If $p\in U$ and $(U,\phi)$ is a smooth chart of $M$ such that $p\in U$, every $A_{p}\in T^{(k,l)}(T_{p}M)$ can be written uniquely as
\[
A_{p}=A^{i_{1}\dots i_{k}}_{j_{1}\dots j_{l}}
\boldsymbol{\partial}_{i_{1}}\biggr|_{p}\otimes\dots\otimes
\boldsymbol{\partial}_{i_{k}}\biggr|_{p}\otimes
\mathbf{d}x^{j_{1}}|_{p}\otimes\dots\otimes
\mathbf{d}x^{j_{l}}|_{p},
\]
where $\phi=(x^{1},\dots,x^{n})$ and each
$A^{i_{1}\dots i_{k}}_{j_{1}\dots j_{l}}\in\mathbb R$ is a component
of \(A_p\) with respect to the tensor frame induced by the chart. When
\(A_p\) varies as a section \(\mathbf A\) over \(U\), these
components become functions
\(A^{i_{1}\dots i_{k}}_{j_{1}\dots j_{l}}\colon U\to\mathbb R\).

 Define $\pi\colon T^{(k,l)}(TM)\longrightarrow M$ by $\pi(A_{p})=p$.

 Take an open cover $\{U_{\alpha}\}_{\alpha\in J}$ by smooth coordinate chart domains for $M$. For each $\alpha\in J$, define $\Phi_{\alpha}\colon \pi^{-1}(U_{\alpha})\longrightarrow U_{\alpha}\times\mathbb{R}^{n^{k+l}}$ by $\Phi_{\alpha}(A_{p})=(p,A^{i_{1}\dots i_{k}}_{j_{1}\dots j_{l}})$. The map $\Phi_{\alpha}$ is bijective because the function $\widetilde{\Phi}_{\alpha}\colon U_{\alpha}\times\mathbb{R}^{n^{k+l}}\longrightarrow \pi^{-1}(U_{\alpha})$ given by \[\widetilde{\Phi}_{\alpha}(p,A^{i_{1}\dots i_{k}}_{j_{1}\dots j_{l}})=A_{p}\] is the inverse of $\Phi_{\alpha}$.

 For each $p\in U_{\alpha}$, $\Phi_{\alpha}|_{T^{(k,l)}(T_{p}M)}\colon T^{(k,l)}(T_{p}M)\longrightarrow \{p\}\times \mathbb{R}^{n^{k+l}}$ is bijective because \[\Phi_{\alpha}|_{T^{(k,l)}(T_{p}M)}(T^{(k,l)}(T_{p}M))=\{p\}\times\mathbb{R}^{n^{k+l}}\] and it is the restriction of a bijective function to its image. Identify the fiber $\{p\}\times\mathbb R^{n^{k+l}}$ with $\mathbb R^{n^{k+l}}$ by the second projection. It suffices to prove that the second component of $\Phi_\alpha$ is linear:

 Let $a\in \mathbb{R}$ and $A_{p},B_{p}\in T^{(k,l)}(T_{p}M)$, with \[A_{p}=A^{i_{1}\dots i_{k}}_{j_{1}\dots j_{l}}(p)\boldsymbol{\partial}_{i_{1}}\biggr|_{p}\otimes\dots\otimes\boldsymbol{\partial}_{i_{k}}\biggr|_{p}\otimes \mathbf{d}x^{j_{1}}|_{p}\otimes\dots\otimes \mathbf{d}x^{j_{l}}|_{p}\] and \[B_{p}=B^{i_{1}\dots i_{k}}_{j_{1}\dots j_{l}}(p)\boldsymbol{\partial}_{i_{1}}\biggr|_{p}\otimes\dots\otimes\boldsymbol{\partial}_{i_{k}}\biggr|_{p}\otimes \mathbf{d}x^{j_{1}}|_{p}\otimes\dots\otimes \mathbf{d}x^{j_{l}}|_{p},\] We have
\[
\Phi_{\alpha}(aA_{p}+B_{p})
=
\left(p,
aA^{i_{1}\dots i_{k}}_{j_{1}\dots j_{l}}
+B^{i_{1}\dots i_{k}}_{j_{1}\dots j_{l}}\right).
\]
Consequently,
\[
\operatorname{pr}_2\Phi_{\alpha}(aA_p+B_p)
=
a\,\operatorname{pr}_2\Phi_{\alpha}(A_p)
+\operatorname{pr}_2\Phi_{\alpha}(B_p),
\]
which is the required linearity.

 Therefore, $\Phi_{\alpha}|_{T^{(k,l)}(T_{p}M)}\colon T^{(k,l)}(T_{p}M)\longrightarrow \{p\}\times \mathbb{R}^{n^{k+l}}$ is a vector space isomorphism.

 Finally, if $\alpha,\beta\in J$ are such that $U_{\alpha}\cap U_{\beta}\neq\varnothing$, we must compute the map $\Phi_{\alpha}\circ \Phi_{\beta}^{-1}$.

 Suppose that $\phi_{\alpha}=(x^{1},\dots,x^{n})$ and $\phi_{\beta}=(y^{1},\dots,y^{n})$.

 If $(p,B^{a_{1}\dots a_{k}}_{b_{1}\dots b_{l}})\in (U_{\alpha}\cap U_{\beta})\times\mathbb{R}^{n^{k+l}}$, we have \[(\Phi_{\alpha}\circ \Phi_{\beta}^{-1})(p,B^{a_{1}\dots a_{k}}_{b_{1}\dots b_{l}})=\Phi_{\alpha}\left(B^{a_{1}\dots a_{k}}_{b_{1}\dots b_{l}}\boldsymbol{\partial}_{a_{1}}\biggr|_{p}\otimes\dots\otimes\boldsymbol{\partial}_{a_{k}}\biggr|_{p}\otimes \mathbf{d}y^{b_{1}}|_{p}\otimes\dots\otimes \mathbf{d}y^{b_{l}}|_{p}\right)\] \[=\left(p,B^{a_{1}\dots a_{k}}_{b_{1}\dots b_{l}}\left(\frac{\partial x^{i_{1}}}{\partial y^{a_{1}}}\biggr|_{p}\right)\dots\left(\frac{\partial x^{i_{k}}}{\partial y^{a_{k}}}\biggr|_{p}\right)\left(\frac{\partial y^{b_{1}}}{\partial x^{j_{1}}}\biggr|_{p}\right)\dots\left(\frac{\partial y^{b_{l}}}{\partial x^{j_{l}}}\biggr|_{p}\right)\right)\] where we have used \[B^{a_{1}\dots a_{k}}_{b_{1}\dots b_{l}}\left(\frac{\partial x^{i_{1}}}{\partial y^{a_{1}}}\biggr|_{p}\right)\dots\left(\frac{\partial x^{i_{k}}}{\partial y^{a_{k}}}\biggr|_{p}\right)\left(\frac{\partial y^{b_{1}}}{\partial x^{j_{1}}}\biggr|_{p}\right)\dots\left(\frac{\partial y^{b_{l}}}{\partial x^{j_{l}}}\biggr|_{p}\right)\boldsymbol{\partial}_{i_{1}}\biggr|_{p}\otimes\dots\otimes\boldsymbol{\partial}_{i_{k}}\biggr|_{p}\otimes \mathbf{d}x^{j_{1}}|_{p}\otimes\dots\otimes \mathbf{d}x^{j_{l}}|_{p}\] \[=B^{a_{1}\dots a_{k}}_{b_{1}\dots b_{l}}\boldsymbol{\partial}_{a_{1}}\biggr|_{p}\otimes\dots\otimes\boldsymbol{\partial}_{a_{k}}\biggr|_{p}\otimes \mathbf{d}y^{b_{1}}|_{p}\otimes\dots\otimes \mathbf{d}y^{b_{l}}|_{p}\] which follows from the chain rule.

 Now note that \[\left(p,B^{a_{1}\dots a_{k}}_{b_{1}\dots b_{l}}\left(\frac{\partial x^{i_{1}}}{\partial y^{a_{1}}}\biggr|_{p}\right)\dots\left(\frac{\partial x^{i_{k}}}{\partial y^{a_{k}}}\biggr|_{p}\right)\left(\frac{\partial y^{b_{1}}}{\partial x^{j_{1}}}\biggr|_{p}\right)\dots\left(\frac{\partial y^{b_{l}}}{\partial x^{j_{l}}}\biggr|_{p}\right)\right)=(p,\tau_{\alpha\beta}(p)B^{a_{1}\dots a_{k}}_{b_{1}\dots b_{l}})\] where \[\tau_{\alpha\beta}(p)=\left(\left(\frac{\partial x^{i_{1}}}{\partial y^{a_{1}}}\biggr|_{p}\right)\dots\left(\frac{\partial x^{i_{k}}}{\partial y^{a_{k}}}\biggr|_{p}\right)\left(\frac{\partial y^{b_{1}}}{\partial x^{j_{1}}}\biggr|_{p}\right)\dots\left(\frac{\partial y^{b_{l}}}{\partial x^{j_{l}}}\biggr|_{p}\right)\right)\in GL(n^{k+l},\mathbb{R})\] is smooth because its component functions are smooth.

 By Lemma~\ref{lema del haz vectorial suave}, there exist a unique topology and a unique smooth structure making $\pi\colon T^{(k,l)}(TM)\longrightarrow M$ a smooth vector bundle of rank $n^{k+l}$ with $\Phi_{\alpha}$ as smooth local trivializations for each $\alpha\in J$.
 The transition formula also holds when comparing a chart from
 this cover with any other smooth chart of $M$. Thus all
 component trivializations induced by smooth charts are smooth,
 and this property characterizes the constructed structure.
 \end{proof}
 \begin{corollary}\label{cor:nociones-fundamentales-carta-suave-cumple-carta-suave}
 Let $M$ be a smooth manifold with or without boundary of dimension $n$. If
 $(U,\phi)$ is a smooth chart of $M$ with
 $\phi=(x^{1},\dots,x^{n})$, then
 $\widetilde{\phi}\colon \pi^{-1}(U)\longrightarrow
 \phi(U)\times\mathbb{R}^{n^{k+l}}$, given by
 \[\widetilde{\phi}\left(A^{i_{1}\dots i_{k}}_{j_{1}\dots j_{l}}\boldsymbol{\partial}_{i_{1}}\biggr|_{p}\otimes\dots\otimes\boldsymbol{\partial}_{i_{k}}\biggr|_{p}\otimes \mathbf{d}x^{j_{1}}|_{p}\otimes\dots\otimes \mathbf{d}x^{j_{l}}|_{p}\right)=(x^{1}(p),\dots,x^{n}(p),A^{i_{1}\dots i_{k}}_{j_{1}\dots j_{l}})\]
 makes $(\pi^{-1}(U),\widetilde{\phi})$ a smooth chart for
 $T^{(k,l)}(TM)$.
 \end{corollary}
 \begin{proof}
 By Proposition~\ref{estructura del haz tensorial}, $\Phi\colon \pi^{-1}(U)\longrightarrow U\times\mathbb{R}^{n^{k+l}}$, given by \[\Phi\left(A^{i_{1}\dots i_{k}}_{j_{1}\dots j_{l}}\boldsymbol{\partial}_{i_{1}}\biggr|_{p}\otimes\dots\otimes\boldsymbol{\partial}_{i_{k}}\biggr|_{p}\otimes \mathbf{d}x^{j_{1}}|_{p}\otimes\dots\otimes \mathbf{d}x^{j_{l}}|_{p}\right)=(p,A^{i_{1}\dots i_{k}}_{j_{1}\dots j_{l}})\] is a smooth local trivialization of the tensor bundle $T^{(k,l)}(TM)$. Moreover, Proposition~\ref{marcos asociados con trivializaciones} says that \[\left(\boldsymbol{\partial}_{i_{1}}\otimes\dots\otimes\boldsymbol{\partial}_{i_{k}}\otimes \mathbf{d}x^{j_{1}}\otimes\dots\otimes \mathbf{d}x^{j_{l}}\right)\] is a smooth local frame defined on $U$ and associated with the local trivialization $\Phi$. Thus Corollary~\ref{carta del haz vectorial} shows that the map $\widetilde{\phi}\colon \pi^{-1}(U)\longrightarrow \phi(U)\times\mathbb{R}^{n^{k+l}}$ makes $(\pi^{-1}(U),\widetilde{\phi})$ a smooth chart of $T^{(k,l)}(TM)$.
 \end{proof}
 This construction of the product of $TM$ and $T^{*}M$ extends to the product of two or more vector bundles, as the reader can prove.
 \begin{exercise}\label{estructura del producto tensorial de haces}
Let $M$ be a smooth manifold with or without boundary and let $\mathbf{E},\mathbf{F}\longrightarrow M$ be two smooth vector bundles of ranks $k$ and $l$, respectively. Prove that the fiberwise tensor product
\[
\mathbf{E}\otimes \mathbf{F} \longrightarrow M,
\qquad (\mathbf{E}\otimes \mathbf{F})_p := \mathbf{E}_p \otimes \mathbf{F}_p,
\]
admits a unique topology and a unique smooth structure making it a smooth vector bundle of rank $kl$ over $M$ for which the tensor trivializations induced by smooth local frames of $\mathbf E$ and $\mathbf F$ are smooth.
\end{exercise}
There is a relation between Hom bundles and tensor products.

\begin{proposition}\label{prop:Hom-tensor}
Let $M$ be a smooth manifold with or without boundary and let $\mathbf{E},\mathbf{F}\longrightarrow M$ be smooth vector bundles of ranks $k$ and $l$, respectively. Then there is a natural vector bundle isomorphism over $M$
\[
\operatorname{Hom}(\mathbf{E},\mathbf{F})\cong \mathbf{E}^{*}\otimes \mathbf{F}.
\]
\end{proposition}

\begin{proof}
Define a bundle map $\Phi\colon \operatorname{Hom}(\mathbf{E},\mathbf{F})\longrightarrow \mathbf{E}^{*}\otimes \mathbf{F}$ fiberwise by the canonical identification
\[
\Phi_p\colon \operatorname{Hom}(\mathbf{E}_p,\mathbf{F}_p)\longrightarrow \mathbf{E}_p^{*}\otimes \mathbf{F}_p,
\qquad
\Phi_p(T_p)=\displaystyle\sum_{i=1}^{k}e^{i}\otimes T_p(e_i),
\]
where $(e_i)_{i=1}^{k}$ is a basis for $\mathbf{E}_p$ and $(e^{i})_{i=1}^{k}$ is the dual basis in $\mathbf{E}_p^{*}$. This definition is independent of the chosen basis, is linear on each fiber, and covers the identity of $M$.

To verify that $\Phi$ is smooth, work in local trivializations. Let $U_{\alpha}\subset M$ be an open set with smooth local frames $(\boldsymbol{\sigma}_{\alpha,i})_{i=1}^{k}$ of $\mathbf{E}$ and $(\boldsymbol{\tau}_{\alpha,j})_{j=1}^{l}$ of $\mathbf{F}$; denote the dual frame of $\mathbf{E}^{*}$ on $U_{\alpha}$ by $(\boldsymbol{\sigma}_{\alpha}^{i})_{i=1}^{k}$.

In these trivializations, every $T_p\in\operatorname{Hom}(\mathbf{E}_p,\mathbf{F}_p)$ is determined by an associated matrix
\[
T_{\alpha}(p)=(T_{\alpha}(p)^{j}_{i})\in M_{l\times k}(\mathbb{R}),
\qquad
T_p(\boldsymbol{\sigma}_{\alpha,i}(p))=\displaystyle\sum_{j=1}^{l}T_{\alpha}(p)^{j}_{i}\boldsymbol{\tau}_{\alpha,j}(p).
\]
The trivialization of $\operatorname{Hom}(\mathbf{E},\mathbf{F})$ on $U_{\alpha}$ is given by
\[
\Theta^{\operatorname{Hom}}_{\alpha}(T_p)=(p,T_{\alpha}(p)).
\]
On the other hand, in $\mathbf{E}^{*}\otimes \mathbf{F}$, any element of the fiber $(\mathbf{E}^{*}\otimes \mathbf{F})_p$ can be written as
\[
\displaystyle\sum_{i=1}^{k}\displaystyle\sum_{j=1}^{l}S^{j}_{i}\boldsymbol{\sigma}_{\alpha}^{i}(p)\otimes \boldsymbol{\tau}_{\alpha,j}(p),
\]
and the associated trivialization is
\[
\Theta^{\mathbf{E}^{*}\otimes \mathbf{F}}_{\alpha}\Bigl(\displaystyle\sum_{\substack{1\leq i\leq k\\1\leq j\leq l}}S^{j}_{i}\boldsymbol{\sigma}_{\alpha}^{i}(p)\otimes \boldsymbol{\tau}_{\alpha,j}(p)\Bigr)
=(p,(S^{j}_{i})).
\]

Now apply $\Phi_p$ to $T_p$. By definition,
\[
\Phi_p(T_p)=\displaystyle\sum_{i=1}^{k}\boldsymbol{\sigma}_{\alpha}^{i}(p)\otimes T_p(\boldsymbol{\sigma}_{\alpha,i}(p))
=\displaystyle\sum_{i=1}^{k}\boldsymbol{\sigma}_{\alpha}^{i}(p)\otimes \Bigl(\displaystyle\sum_{j=1}^{l}T_{\alpha}(p)^{j}_{i}\boldsymbol{\tau}_{\alpha,j}(p)\Bigr).
\]
Reordering the sum,
\[
\Phi_p(T_p)=\displaystyle\sum_{i=1}^{k}\displaystyle\sum_{j=1}^{l}T_{\alpha}(p)^{j}_{i}\boldsymbol{\sigma}_{\alpha}^{i}(p)\otimes \boldsymbol{\tau}_{\alpha,j}(p).
\]

In the trivialization of $\mathbf{E}^{*}\otimes \mathbf{F}$, this corresponds exactly to
\[
\Theta^{\mathbf{E}^{*}\otimes \mathbf{F}}_{\alpha}\bigl(\Phi_p(T_p)\bigr)=(p,(T_{\alpha}(p)^{j}_{i})).
\]

In other words, passing through the two local trivializations, the map
\[
\Theta^{\mathbf{E}^{*}\otimes \mathbf{F}}_{\alpha}\circ \Phi \circ \bigl(\Theta^{\operatorname{Hom}}_{\alpha}\bigr)^{-1}\colon
U_{\alpha}\times M_{l\times k}(\mathbb{R})\to U_{\alpha}\times M_{l\times k}(\mathbb{R})
\]
is simply given by
\[
(p,(T^{j}_{i}))\longmapsto (p,(T^{j}_{i})),
\]
which is the identity, showing that $\Phi$ is smooth.

Thus $\Phi$ is smooth and locally an isomorphism of trivializations. Since each $\Phi_p$ is linear and bijective, $\Phi$ is a bijective bundle homomorphism. By Exercise \ref{homomorfismo biyectivo es isomorfismo}, we conclude that
\[
\Phi:\operatorname{Hom}(\mathbf{E},\mathbf{F})\xrightarrow{\cong}\mathbf{E}^{*}\otimes \mathbf{F}
\]
is a smooth bundle isomorphism of rank $kl$.

The naturality of the construction follows directly from the fact that everything was defined fiberwise using the natural isomorphism $\operatorname{Hom}(V,W)\cong V^{*}\otimes W$.
\end{proof}

Sections of tensor bundles are tensor fields. They describe, among other objects, metrics, differential forms, curvature, and tensors arising in physical models.
 \begin{definition}\label{def:nociones-fundamentales-campo-tensorial-de-tipo}\index{tensor field!mixed}
 Let $M$ be a smooth manifold with or without boundary.
 A section
 \[
 \mathbf{A}\in \Gamma\bigl(T^{(k,l)}(TM)\bigr)
 \]
 is called a \textbf{tensor field of type $(k,l)$} on $M$.

 In particular:
 \begin{itemize}
 \item If $k=0$, we say that $\mathbf{A}$ is a \textbf{$l$-covariant tensor field}.
 \item If $l=0$, we say that $\mathbf{A}$ is a \textbf{$k$-contravariant tensor field}.
 \end{itemize}
\end{definition}
\begin{remark}\label{obs:nociones-fundamentales-campos-vectoriales}
 We use the following classical notation:
 \begin{itemize}
 \item We denote the space of \textbf{vector fields} on $M$ by $\mathfrak{X}(M):=\Gamma(T^{(1,0)}(TM))$.
 \item We denote the space of \textbf{1-forms} on $M$ by $\Omega^{1}(M):=\Gamma(T^{(0,1)}(TM))$.
 \end{itemize}
 The notation $\Omega^{1}(M)$ will be explained in greater detail below.
\end{remark}
 The following criteria determine when a tensor field is smooth and when a map is induced by a smooth tensor field.

 \begin{proposition}[Smoothness criterion for tensor fields]\label{criterio de suavidad para campos tensoriales}\index{smoothness criterion for tensor fields}
 Let $M$ be a smooth manifold with or without boundary, and let
 $\mathbf{F}\colon M\longrightarrow T^{(k,l)}(TM)$ be a rough $(k,l)$-tensor field. Then $\mathbf{F}$ is smooth on an open set
 $U\subseteq M$ if and only if, for any smooth covector fields
 $\boldsymbol{\omega}^{1},\dots,\boldsymbol{\omega}^{k}$ and smooth vector fields $\mathbf{X}_{1},\dots,\mathbf{X}_{l}$ defined on
 $U$, the real-valued function
 \[
 \mathbf{F}(\boldsymbol{\omega}^{1},\dots,\boldsymbol{\omega}^{k},
 \mathbf{X}_{1},\dots,\mathbf{X}_{l})(p)
 =\mathbf{F}_{p}(\boldsymbol{\omega}^{1}|_{p},\dots,
 \boldsymbol{\omega}^{k}|_{p},\mathbf{X}_{1}|_{p},\dots,
 \mathbf{X}_{l}|_{p})
 \]
 is smooth on $U$.
 \end{proposition}
\begin{proof}
If \(\mathbf F\) is smooth, its evaluation on smooth fields is a smooth function, since in coordinates it is a finite sum of products of the smooth components of \(\mathbf F\) and of the fields on which it is evaluated.

Conversely, suppose that every indicated evaluation is smooth, and fix \(p\in U\). Choose a chart \((V,x^1,\dots,x^n)\) with \(p\in V\subseteq U\), and an open set \(W\) such that \(p\in W\Subset V\). Take \(\rho\in C_c^\infty(V)\) equal to one on \(W\). The fields
\[
\rho\,\mathbf d x^1,\dots,\rho\,\mathbf d x^n,
\qquad
\rho\,\boldsymbol{\partial}_1,\dots,\rho\,\boldsymbol{\partial}_n
\]
extend by zero to smooth fields defined on \(U\). For indices \(i_1,\dots,i_k,j_1,\dots,j_l\), the hypothesis implies that
\[
\mathbf F\bigl(
\rho\,\mathbf d x^{i_1},\dots,\rho\,\mathbf d x^{i_k},
\rho\,\boldsymbol{\partial}_{j_1},\dots,
\rho\,\boldsymbol{\partial}_{j_l}
\bigr)
\]
is smooth on \(U\). On \(W\), where \(\rho=1\), this function is precisely the coordinate component
\[
F^{\,i_1\dots i_k}_{j_1\dots j_l}
:=
\mathbf F\bigl(
\mathbf d x^{i_1},\dots,\mathbf d x^{i_k},
\boldsymbol{\partial}_{j_1},\dots,\boldsymbol{\partial}_{j_l}
\bigr).
\]
Thus all components of \(\mathbf F\) are smooth in a neighborhood of \(p\). Since \(p\) was arbitrary, \(\mathbf F\) is smooth on \(U\).
\end{proof}

 \begin{lemma}[Tensor characterization lemma]\label{lema de caracterización tensorial}\index{tensor characterization lemma}
 Let $M$ be a smooth manifold with or without boundary.

 A map
 \[
 F:\underbrace{\Omega^{1}(M)\times\dots\times\Omega^{1}(M)}_{k\text{ times}}
 \times
 \underbrace{\mathfrak{X}(M)\times\dots\times \mathfrak{X}(M)}_{l\text{ times}}
 \longrightarrow C^{\infty}(M)
 \]
 is induced by a smooth $(k,l)$-tensor field if and only if it is multilinear over $C^{\infty}(M)$.

 For the vector-valued statement, a map
 \[
 F:\underbrace{\Omega^{1}(M)\times\dots\times\Omega^{1}(M)}_{k\text{ times}}
 \times
 \underbrace{\mathfrak{X}(M)\times\dots\times \mathfrak{X}(M)}_{l\text{ times}}
 \longrightarrow \mathfrak{X}(M)
 \]
 is induced by a smooth $(k+1,l)$-tensor field if and only if it is multilinear over $C^{\infty}(M)$.
\end{lemma}

\begin{proof}
We first prove a locality observation valid for smooth sections of any finite-rank vector bundle. Let $\mathcal F$ be a map that is
$C^\infty(M)$-linear in one of its arguments, and let $\mathbf{s}$ be a section that vanishes in a neighborhood of $p$. By
Proposition~\ref{funciones flan}, there exists $\vartheta\in C^\infty(M)$ with
$\vartheta(p)=1$ and support contained in that neighborhood. Then
\[
 \mathcal F(\ldots,\mathbf{s},\ldots)(p)
 =\vartheta(p)\mathcal F(\ldots,\mathbf{s},\ldots)(p)
 =\mathcal F(\ldots,\vartheta \mathbf{s},\ldots)(p)=0.
\]
If two sections $\mathbf{s}$ and $\widetilde{\mathbf{s}}$ have the same value at $p$, choose a trivialization around $p$ and a cutoff function $\rho$ that equals one in a neighborhood of $p$ and has support within the trivialization. The section $\mathbf{r}:=\mathbf{s}-\widetilde{\mathbf{s}}$ satisfies
$\mathbf{r}(p)=0$. The term $(1-\rho)\mathbf{r}$ vanishes near $p$, whereas, after extending the cutoff coefficients and the elements of a local frame by zero, we can write
\[
 \rho \mathbf{r}=\displaystyle\sum_{i=1}^{m}a_i \mathbf{e}_i,
 \qquad a_i\in C^\infty(M),\qquad a_i(p)=0.
\]
The preceding locality observation and $C^\infty(M)$-linearity give
$\mathcal F(\ldots,\mathbf{r},\ldots)(p)=0$. Thus the value at $p$ depends only on the values at $p$ of all the arguments.

Now suppose that
\[
 F:\Omega^1(M)^k\times\mathfrak X(M)^l\longrightarrow C^\infty(M)
\]
is $C^\infty(M)$-multilinear. Given
$\alpha^1,\ldots,\alpha^k\in T_p^*M$ and
$v_1,\ldots,v_l\in T_pM$, choose smooth extensions and define
\[
 \mathbf{T}_p(\alpha^1,\ldots,\alpha^k,v_1,\ldots,v_l)
 :=F(\widetilde{\boldsymbol{\alpha}}^1,\ldots,
       \widetilde{\boldsymbol{\alpha}}^k,
       \widetilde{\mathbf{v}}_1,\ldots,\widetilde{\mathbf{v}}_l)(p).
\]
The preceding observation shows that this definition is independent of the extensions; the multilinearity of $\mathbf{T}_p$ is inherited from that of $F$. In a chart, choose global extensions of the coordinate coframe and frame that agree with them on a smaller open set. The components of
$\mathbf{T}$ on that open set are the restrictions of the smooth functions obtained by evaluating $F$ on those extensions. Hence $p\mapsto \mathbf{T}_p$ is a smooth
$(k,l)$-tensor field and, by construction, induces $F$. Uniqueness follows by evaluating the field on arbitrary fiber values. Conversely, evaluating a smooth tensor field on smooth sections produces a smooth function and is $C^\infty(M)$-multilinear.

If the codomain of $F$ is $\mathfrak X(M)$, we instead define
\[
 \mathbf{T}_p(\alpha^1,\ldots,\alpha^k,v_1,\ldots,v_l)
 :=F(\widetilde{\boldsymbol{\alpha}}^1,\ldots,
       \widetilde{\boldsymbol{\alpha}}^k,
       \widetilde{\mathbf{v}}_1,\ldots,\widetilde{\mathbf{v}}_l)(p)\in T_pM.
\]
The same locality argument gives independence of the extensions. The components of the vector field on the right-hand side are smooth; evaluating on local coframes and frames shows that $\mathbf{T}$ is a smooth
$(k+1,l)$-tensor field. The converse implication and uniqueness follow by evaluation, exactly as in the scalar case.
\end{proof}

 We now define the pullback of a tensor field.
 \begin{definition}\label{def:nociones-fundamentales-pullback-de-bajo-en}\index{pullback!of a covariant tensor}
 Let $M$ and $N$ be smooth manifolds with or without boundary, and suppose that
 $F\colon M\longrightarrow N$ is a smooth map.
 For each point $p\in M$ and each $k$-covariant tensor
 \(\alpha\in T^{(0,k)}(T_{F(p)}N)\),
 we define a tensor
 \[
 dF_{p}^{*}(\alpha)\in T^{(0,k)}(T_{p}M),
 \]
 called the \textbf{pullback of $\alpha$ by $F$ at $p$}, by
 \[
 dF_{p}^{*}(\alpha)(v_{1},\dots,v_{k})
 =\alpha\bigl(dF_{p}(v_{1}),\dots,dF_{p}(v_{k})\bigr),
 \]
 for any $v_{1},\dots,v_{k}\in T_{p}M$.

 If $\mathbf{A}$ is a $k$-covariant tensor field on $N$, that is,
 \(\mathbf{A}\in \Gamma(T^{(0,k)}(TN))\),
 we define the tensor field on $M$
 \[
 F^{*}\mathbf{A} \in \Gamma\bigl(T^{(0,k)}(TM)\bigr),
 \]
 called the \textbf{pullback of $\mathbf{A}$ by $F$}, by
 \[
 (F^{*}\mathbf{A})_{p} = dF_{p}^{*}(\mathbf{A}_{F(p)}),
 \]
 that is,
 \[
 (F^{*}\mathbf{A})_{p}(v_{1},\dots,v_{k})
 = \mathbf{A}_{F(p)}\bigl(dF_{p}(v_{1}),\dots,dF_{p}(v_{k})\bigr),
 \]
 for each $v_{1},\dots,v_{k}\in T_{p}M$.
\end{definition}

 The pullback has the following properties:
 \begin{proposition}\label{propiedades del pullback}
 Let $M$, $N$, and $P$ be smooth manifolds with or without boundary. Suppose that
 $F\colon M\longrightarrow N$ and $G\colon N\longrightarrow P$ are smooth maps. Let \(k,l,r\in\mathbb N_0\),
 \(\mathbf A,\mathbf B\in\Gamma(T^{(0,k)}(TN))\),
 \(\mathbf C\in\Gamma(T^{(0,l)}(TN))\),
 \(\mathbf D\in\Gamma(T^{(0,r)}(TP))\), and \(f\in C^\infty(N)\).
 \begin{enumerate}[label=(\alph*)]
 \item $F^{*}(f\mathbf{A})=(f\circ F)F^{*}\mathbf{A}$.
 \item $F^{*}(\mathbf{A}\otimes \mathbf{C})=F^{*}\mathbf{A}\otimes F^{*}\mathbf{C}$.
 \item $F^{*}(\mathbf{A}+\mathbf{B})=F^{*}\mathbf{A}+F^{*}\mathbf{B}$.
 \item $F^{*}\mathbf{A}$ is a smooth tensor field.
 \item $(G\circ F)^{*}\mathbf{D}=F^{*}(G^{*}\mathbf{D})$.
 \item $(\operatorname{Id}_{N})^{*}\mathbf{A}=\mathbf{A}$.
 \end{enumerate}

 \end{proposition}
\begin{proof}
The first equality follows by evaluating on \(p\in M\):
\[
\bigl(F^*(f\mathbf A)\bigr)_p(v_1,\dots,v_k)
=f(F(p))\,\mathbf A_{F(p)}(dF_pv_1,\dots,dF_pv_k).
\]
The definition of the tensor product proves (b), and the linearity of
\(dF_p^*\) proves (c). For (d), if in coordinates
\[
\mathbf A=A_{i_1\dots i_k}\,
\mathbf d y^{i_1}\otimes\cdots\otimes\mathbf d y^{i_k},
\]
then
\[
F^*\mathbf A=(A_{i_1\dots i_k}\circ F)\,
d(y^{i_1}\circ F)\otimes\cdots\otimes d(y^{i_k}\circ F),
\]
whose coefficients are smooth.

Substituting the chain rule \(d(G\circ F)_p=dG_{F(p)}\circ dF_p\) into the pointwise definition gives
\((G\circ F)^*\mathbf D=F^*(G^*\mathbf D)\). Finally,
\(d(\operatorname{Id}_N)_q=\operatorname{Id}_{T_qN}\), which proves (f).
\end{proof}

 A related concept, the pushforward, requires $F$ to be a diffeomorphism as well.

 To introduce it, we first define the pushforward for vector fields:
 \begin{definition}\label{def:nociones-fundamentales-pushforward}\index{pushforward}
 Let $M$ and $N$ be smooth manifolds with or without boundary, and let $F\colon M\longrightarrow N$ be a diffeomorphism. For each $\mathbf{X}\in\mathfrak{X}(M)$, we define $F_{*}\mathbf X:=dF\circ \mathbf{X}\circ F^{-1}$ and call it the \textbf{pushforward} of $\mathbf{X}$ by $F$. Pointwise, we have \[(F_{*}\mathbf{X})_{q}=dF_{F^{-1}(q)}(\mathbf{X}_{F^{-1}(q)})\] for each $q\in N$.
 \end{definition}
 The most general form of pullback and pushforward is given by the following proposition:
 \begin{proposition}\label{pullback y pushforward}
 Let $M$, $N$, and $P$ be smooth manifolds with or without boundary, and let
 \(F\colon M\longrightarrow N\) and \(G\colon N\longrightarrow P\)
 be diffeomorphisms. For each \(k,l\in\mathbb N_0\), there are isomorphisms associated with \(F\), called the \textbf{pushforward} and \textbf{pullback},
 \[
 F_{*}\colon \Gamma(T^{(k,l)}(TM))\longrightarrow \Gamma(T^{(k,l)}(TN)),
 \qquad
 F^{*}\colon \Gamma(T^{(k,l)}(TN))\longrightarrow \Gamma(T^{(k,l)}(TM)),
 \]
 such that $F^{*}$ agrees with the pullback of covariant tensor fields,
 $F_{*}$ agrees with the pushforward of vector fields,
 and the following conditions hold:
 \begin{enumerate}[label=(\alph*)]
 \item $F_{*}=(F^{*})^{-1}$.
 \item $F^{*}(\mathbf{A}\otimes\mathbf B)=F^{*}\mathbf{A}\otimes F^{*}\mathbf B$ for any tensor fields \(\mathbf A\) and \(\mathbf B\) on \(N\).
 \item $(G\circ F)_{*}=G_{*}\circ F_{*}$.
 \item $(G\circ F)^{*}=F^{*}\circ G^{*}$.
 \item $(\text{Id}_{M})^{*}=(\text{Id}_{M})_{*}=\text{Id}\colon \Gamma(T^{(k,l)}(TM))\longrightarrow \Gamma(T^{(k,l)}(TM))$.
 \item $F^{*}\bigl(\mathbf{A}(\mathbf{X}_{1},\dots, \mathbf{X}_{k})\bigr)=(F^{*}\mathbf{A})\bigl((F^{-1})_{*}\mathbf X_{1},\dots,(F^{-1})_{*}\mathbf X_{k}\bigr)$ for each $\mathbf{A}\in \Gamma(T^{(0,k)}(TN))$ and any $\mathbf{X}_{1},\dots,\mathbf{X}_{k}\in \mathfrak{X}(N)$.
 \end{enumerate}
\end{proposition}
\begin{proof}
Let \(q\in N\), set \(p=F^{-1}(q)\), and let
\(\mathbf T\in\Gamma(T^{(k,l)}(TM))\). For
\(\beta^1,\dots,\beta^k\in T_q^*N\) and
\(w_1,\dots,w_l\in T_qN\), define
\[
\begin{aligned}
(F_*\mathbf T)_q(\beta^1,\dots,\beta^k,w_1,\dots,w_l)
:={}&\mathbf T_p\bigl(dF_p^*\beta^1,\dots,dF_p^*\beta^k,\\
&\hspace{27mm}(dF_p)^{-1}w_1,\dots,(dF_p)^{-1}w_l\bigr).
\end{aligned}
\]
In charts, its components are finite sums of products of components of
\(\mathbf T\), \(dF\), and \(dF^{-1}\); it therefore defines a smooth field.
For a vector field, it agrees with the pushforward already defined. Set
\(F^*:=(F^{-1})_*\); on covariant tensors, this is the usual definition of the pullback.

Substituting the identity \(d(F^{-1})_q=(dF_p)^{-1}\) into the formula shows that \(F_*\) and \(F^*\) are inverses. The same formula and the definition of the tensor product prove (b). The chain rule and
\[
\bigl(d(G\circ F)_p\bigr)^{-1}
=(dF_p)^{-1}\circ(dG_{F(p)})^{-1}
\]
give \((G\circ F)_*=G_*\circ F_*\). Applying this equality to the inverses yields
\[
(G\circ F)^*=(F^{-1}\circ G^{-1})_*=F^*\circ G^*.
\]
The identity case is immediate from the formula.

Finally, if \(p\in M\), the right-hand side of (f), evaluated on \(p\), is
\[
\mathbf A_{F(p)}
\bigl(dF_p(dF_p)^{-1}\mathbf X_1(F(p)),\dots,
      dF_p(dF_p)^{-1}\mathbf X_k(F(p))\bigr),
\]
which agrees with
\(\mathbf A(\mathbf X_1,\dots,\mathbf X_k)(F(p))\), that is, with the left-hand side evaluated on \(p\).
\end{proof}
\begin{definition}[Interior multiplication]\label{def:multiplicacion-interior}\index{interior multiplication}
Let $M$ be a smooth manifold with or without boundary, and let $k,l\geq 0$ be nonnegative integers.

\begin{enumerate}[label=(\alph*)]
 \item Let $\mathbf{T}_{p}\in T^{(k,l+1)}(T_{p}M)$ and $\mathbf{V}\in T_{p}M$.
 We define the \emph{interior multiplication} of $\mathbf{V}$ into $\mathbf{T}_{p}$ to be the tensor
 \[
 \iota_{\mathbf{V}}\mathbf{T}_{p} \in T^{(k,l)}(T_{p}M),
 \]
 given by
 \[
 (\iota_{\mathbf{V}}\mathbf{T}_{p})(\omega^{1},\dots,\omega^{k},X_{1},\dots,X_{l})
 :=\mathbf{T}_{p}(\omega^{1},\dots,\omega^{k},\mathbf{V},X_{1},\dots,X_{l}),
 \]
 for all $\omega^{1},\dots,\omega^{k}\in \mathbf{T}_{p}^{*}M$ and $X_{1},\dots,X_{l}\in T_{p}M$.

 \item If $\mathbf{T}\in\Gamma\big(T^{(k,l+1)}(TM)\big)$ and $\mathbf{V}\in\mathfrak{X}(M)$, we define the section
 \[
 \iota_{\mathbf{V}}\mathbf{T} \in\Gamma\big(T^{(k,l)}(TM)\big),
 \]
 by
 \[
 (\iota_{\mathbf{V}}\mathbf{T})(p):=\iota_{\mathbf{V}(p)}\big(\mathbf{T}(p)\big),\qquad p\in M.
 \]
\end{enumerate}
\end{definition}
\begin{proposition}[Interior multiplication defines a tensor field]\label{prop:multiplicacion-interior-campo}\index{interior multiplication defines a tensor field}
Let $M$ be a smooth manifold with or without boundary, $\mathbf{T}\in\Gamma\big(T^{(k,l+1)}(TM)\big)$, and $\mathbf{V}\in\mathfrak{X}(M)$. Then the assignment
\[
(\iota_{\mathbf{V}}\mathbf{T})(\omega^{1},\dots,\omega^{k},X_{1},\dots,X_{l})
:=\mathbf{T}(\omega^{1},\dots,\omega^{k},\mathbf{V},X_{1},\dots,X_{l})
\]
defines a smooth $(k,l)$-tensor field $\iota_{\mathbf{V}}\mathbf{T}\in\Gamma\big(T^{(k,l)}(TM)\big)$.
\end{proposition}

\begin{proof}
Consider the map
\[
F_{\mathbf{V},\mathbf{T}}:
\underbrace{\Omega^{1}(M)\times\dots\times\Omega^{1}(M)}_{k\text{ times}}
\times
\underbrace{\mathfrak{X}(M)\times\dots\times\mathfrak{X}(M)}_{l\text{ times}}
\longrightarrow C^{\infty}(M).
\]
defined by
\[
F_{\mathbf{V},\mathbf{T}}(\omega^{1},\dots,\omega^{k},X_{1},\dots,X_{l})
= \mathbf{T}(\omega^{1},\dots,\omega^{k},\mathbf{V},X_{1},\dots,X_{l}).
\]

Since $\mathbf{T}$ is a smooth $(k,l+1)$-tensor field, the map
\[
(\omega^{1},\dots,\omega^{k},Y,X_{1},\dots,X_{l})\longmapsto
\mathbf{T}(\omega^{1},\dots,\omega^{k},Y,X_{1},\dots,X_{l})
\]
is multilinear over $C^{\infty}(M)$ in its arguments. Fixing $Y=\mathbf{V}$, it follows that $F_{\mathbf{V},\mathbf{T}}$ is multilinear over $C^{\infty}(M)$.

By Lemma~\ref{lema de caracterización tensorial}, there exists a smooth $(k,l)$-tensor field $\mathbf{S}\in\Gamma\big(T^{(k,l)}(TM)\big)$ such that
\[
\mathbf{S}(\omega^{1},\dots,\omega^{k},X_{1},\dots,X_{l})
=F_{\mathbf{V},\mathbf{T}}(\omega^{1},\dots,\omega^{k},X_{1},\dots,X_{l}).
\]

Pointwise, we have
\[
\mathbf{S}(p)(\omega^{1}_{p},\dots,\omega^{k}_{p},X_{1,p},\dots,X_{l,p})
=\mathbf{T}_{p}(\omega^{1}_{p},\dots,\omega^{k}_{p},\mathbf{V}_{p},X_{1,p},\dots,X_{l,p})
=(\iota_{\mathbf{V}(p)}\mathbf{T}_{p})(\omega^{1}_{p},\dots,X_{l,p}),
\]
so $\mathbf{S}=\iota_{\mathbf{V}}\mathbf{T}$ as sections. Therefore $\iota_{\mathbf{V}}\mathbf{T}$ is a smooth $(k,l)$-tensor field.
\end{proof}
\begin{remark}\label{expresion local multiplicacion interior}
 If $\mathbf{T}\in\Gamma\big(T^{(k,l+1)}(TM)\big)$ and $\mathbf{V}\in\mathfrak{X}(M)$ with $\mathbf{V}=V^{a}\mathbf{E}_{a}$, then the interior multiplication $\iota_{\mathbf{V}}\mathbf{T}$ has local components
\[
(\iota_{\mathbf{V}}\mathbf{T})^{i_{1}\dots i_{k}}_{j_{1}\dots j_{l}}
=(\iota_{\mathbf{V}}\mathbf{T})\big(\boldsymbol{\varepsilon}^{i_{1}},\dots,\boldsymbol{\varepsilon}^{i_{k}},\mathbf{E}_{j_{1}},\dots,\mathbf{E}_{j_{l}}\big)
\]
\[
=\mathbf{T}\big(\boldsymbol{\varepsilon}^{i_{1}},\dots,\boldsymbol{\varepsilon}^{i_{k}},\mathbf{V},\mathbf{E}_{j_{1}},\dots,\mathbf{E}_{j_{l}}\big)
\]
\[
=\mathbf{T}\big(\boldsymbol{\varepsilon}^{i_{1}},\dots,\boldsymbol{\varepsilon}^{i_{k}},V^{a}\mathbf{E}_{a},\mathbf{E}_{j_{1}},\dots,\mathbf{E}_{j_{l}}\big)
\]
\[
=V^{a}\mathbf{T}\big(\boldsymbol{\varepsilon}^{i_{1}},\dots,\boldsymbol{\varepsilon}^{i_{k}},\mathbf{E}_{a},\mathbf{E}_{j_{1}},\dots,\mathbf{E}_{j_{l}}\big)
\]
\[
=V^{a}T^{i_{1}\dots i_{k}}_{aj_{1}\dots j_{l}}.
\]
Thus, in any local frame, the operation $\iota_{\mathbf{V}}$ is obtained by contracting the contravariant index of $\mathbf{V}$ with the first covariant index of $\mathbf{T}$.
\end{remark}

 We can introduce two classes of tensors that will be particularly useful and that are used to define Riemannian metrics and differential forms. The first is the class of symmetric tensors:
 \begin{definition}\label{def:nociones-fundamentales-simetrico}\index{symmetric}
 Let $V$ be a finite-dimensional vector space. A $k$-covariant tensor $\alpha$ on $V$ is said to be \textbf{symmetric} if its value is unchanged when any two of its arguments are interchanged:
 \[\alpha(v_{1},\dots,v_{i},\dots,v_{j},\dots,v_{k})=\alpha(v_{1},\dots,v_{j},\dots,v_{i},\dots,v_{k})\] for any $1\leq i<j\leq k$. We denote the space of symmetric $k$-covariant tensors by $\Sigma^{k}(V^{*})\subseteq T^{k}(V^{*})$.
 \end{definition}
 \begin{definition}\label{def:nociones-fundamentales-producto-simetrico}\index{symmetric product}
 Let $V$ be a finite-dimensional vector space. If $\alpha\in \Sigma^{k}(V^{*})$ and $\beta\in \Sigma^{l}(V^{*})$, we define their \textbf{symmetric product} to be the $(k+l)$-covector given by \[\alpha\beta=\frac{1}{(k+l)!}\displaystyle\sum_{\sigma\in S_{k+l}}{}^{\sigma}(\alpha\otimes\beta)\] where $^{\sigma}(\alpha\otimes\beta)(v_{1},\dots,v_{k+l})=\alpha\otimes\beta(v_{\sigma(1)},\dots,v_{\sigma(k+l)})$.
 \end{definition}
 \begin{definition}\label{def:nociones-fundamentales-campo-tensorial-simetrico}\index{symmetric tensor field}
 A \textbf{symmetric tensor field} on a smooth manifold with or without boundary is a covariant tensor field whose value at each point is a symmetric tensor.
 \end{definition}
 The second class consists of alternating tensors:
 \begin{definition}\label{def:nociones-fundamentales-alternante}\index{covariant tensor!alternating}
 If $V$ is a finite-dimensional real vector space, a $k$-covariant tensor $\alpha$ is said to be \textbf{alternating} or \textbf{antisymmetric} if, for any $v_{1},\dots,v_{k}\in V$ and each $i,j\in \{1,\dots,k\}$ with $i\neq j$, we have \[\alpha(v_{1},\dots,v_{i},\dots,v_{j},\dots,v_{k})=-\alpha(v_{1},\dots,v_{j},\dots,v_{i},\dots,v_{k}).\]

 The space of alternating $k$-covariant tensors is denoted by $\Lambda^{k}(V^{*})\subseteq T^{(0,k)}(V)$.
 \end{definition}
 Like the tensor spaces above, the space of alternating $k$-covariant tensors is a finite-dimensional vector space. Alternating $k$-covariant tensors are often called \textit{\textbf{$k$-covectors}}.
 \begin{definition}\label{def:nociones-fundamentales-espacio-vectorial-dimension-finita-base-multiindice}\index{alternating covector!associated with a multi-index}
 Let $V$ be a finite-dimensional vector space, and suppose that $(\varepsilon^{1},\dots,\varepsilon^{n})$ is any basis for $V^{*}$. For each multi-index $I=(i_{1},\dots,i_{k})$ of length $k$ such that $1\leq i_{1},\dots,i_{k}\leq n$, define $\varepsilon^{I}=\varepsilon^{i_{1}\dots i_{k}}$ by \[\varepsilon^{I}(v_{1},\dots,v_{k})=\det \begin{bmatrix}
 \varepsilon^{i_{1}}(v_{1}) & \cdots & \varepsilon^{i_{1}}(v_{k}) \\
 \vdots & \ddots & \vdots \\
 \varepsilon^{i_{k}}(v_{1}) & \cdots & \varepsilon^{i_{k}}(v_{k})
 \end{bmatrix}=\det\begin{bmatrix}
 v_{1}^{i_{1}} & \cdots & v_{k}^{i_{1}} \\
 \vdots & \ddots & \vdots \\
 v_{1}^{i_{k}} & \cdots & v_{k}^{i_{k}}
 \end{bmatrix}\]
 \end{definition}
 \begin{remark}\label{obs:nociones-fundamentales-multiindice}
 For each multi-index $I$, $\varepsilon^{I}\in \Lambda^{k}(V^{*})$.
 \end{remark}
 \begin{proposition}\label{prop:nociones-fundamentales-espacio-vectorial-dimension-finita-base-coleccion}
 Let $V$ be a vector space of finite dimension $n$. If $(\varepsilon^{i})$ is any basis for $V^{*}$, then, for each $k\leq n$, the collection of $k$-covectors \[\mathscr{B}=\{\varepsilon^{I}\mid \text{$I$ is an increasing multi-index of length $k$}\}\] is a basis for $\Lambda^{k}(V^{*})$. Therefore $\dim \Lambda^{k}(V^{*})=\displaystyle{n\choose k}.$
 \end{proposition}
\begin{proof}
Let \((e_1,\dots,e_n)\) be the basis of \(V\) dual to \((\varepsilon^1,\dots,\varepsilon^n)\). If \(I=(i_1<\cdots<i_k)\) and \(J=(j_1<\cdots<j_k)\), the definition by determinants gives
\[
\varepsilon^I(e_{j_1},\dots,e_{j_k})
=
\begin{cases}
1,& I=J,\\
0,& I\neq J.
\end{cases}
\]
Consequently, the elements of \(\mathscr B\) are linearly independent.

Now let \(\alpha\in\Lambda^k(V^*)\) and define
\[
\widetilde\alpha
:=
\displaystyle\sum_{1\leq i_1<\cdots<i_k\leq n}
\alpha(e_{i_1},\dots,e_{i_k})\,
\varepsilon^{i_1\dots i_k}.
\]
If \(j_1,\dots,j_k\) are distinct, there exists a permutation \(\sigma\in S_k\) that arranges them in increasing order. The alternating property of \(\alpha\) and \(\widetilde\alpha\), together with the preceding equality, shows that
\[
\widetilde\alpha(e_{j_1},\dots,e_{j_k})
=
\alpha(e_{j_1},\dots,e_{j_k}).
\]
If two indices coincide, both sides are zero. By multilinearity, \(\widetilde\alpha=\alpha\) on all of \(V^k\), so \(\mathscr B\) spans \(\Lambda^k(V^*)\). Finally, the number of strictly increasing multi-indices of length \(k\) is \(\binom nk\), which proves the dimension formula.
\end{proof}

 We define an operation on $\Lambda^{k}(V^{*})$ as follows:
 \begin{definition}\label{def:nociones-fundamentales-producto-cuna}\index{wedge product} Let $V$ be a finite-dimensional vector space.
 If $\omega\in \Lambda^{k}(V^{*})$ and $\eta\in \Lambda^{l}(V^{*})$, we define their \textbf{wedge product} or \textbf{exterior product} to be the $(k+l)$-covector given by \[\omega\wedge \eta=\frac{1}{k!l!}\displaystyle\sum_{\sigma\in S_{k+l}}(\sign \sigma)(^{\sigma}(\omega\otimes\eta)).\]
 Here the combinatorial factor is included for normalization.
 \end{definition}
 The wedge product has the following properties:
 \begin{proposition}
 \label{propiedades del producto cuña} Suppose that $\omega,\omega',\eta,\eta'$ and $\xi$ are multicovectors on a finite-dimensional vector space $V$.
 \begin{enumerate}[label=(\alph*)]
 \item The wedge product is bilinear over $\mathbb{R}$ and associative.
 \item If $\omega\in \Lambda^{k}(V^{*})$ and $\eta\in\Lambda^{l}(V^{*})$, \[\omega\wedge\eta=(-1)^{kl}\eta\wedge\omega.\]
 \item If $\omega^{1},\dots,\omega^{k}\in V^{*}$ and $v_{1},\dots,v_{k}\in V$, then \[\omega^{1}\wedge\dots\wedge\omega^{k}(v_{1},\dots,v_{k})=\det(\omega^{j}(v_{i})).\]

 \end{enumerate}
 \end{proposition}
\begin{proof}
In the sum defining \(\omega\wedge\eta\), internal permutations of the first \(k\) arguments and the last \(l\) arguments yield equal terms because \(\omega\) and \(\eta\) are alternating. Grouping these \(k!l!\) terms gives the shuffle formula
\[
\begin{aligned}
(\omega\wedge\eta)(v_1,\dots,v_{k+l})
=
\displaystyle\sum_{\sigma\in\operatorname{Sh}(k,l)}
\operatorname{sgn}(\sigma)\,
&\omega(v_{\sigma(1)},\dots,v_{\sigma(k)})\\
{}\cdot&
\eta(v_{\sigma(k+1)},\dots,v_{\sigma(k+l)}),
\end{aligned}
\]
where \(\operatorname{Sh}(k,l)\) is the set of permutations that preserve the relative order within the two blocks. This expression immediately proves bilinearity.

For associativity, let the degrees of \(\omega,\eta,\xi\) be \(k,l,m\), respectively. Applying the shuffle formula twice turns both \((\omega\wedge\eta)\wedge\xi\) and \(\omega\wedge(\eta\wedge\xi)\) into
\[
\displaystyle\sum_{\sigma\in\operatorname{Sh}(k,l,m)}
\operatorname{sgn}(\sigma)\,
\omega(v_{\sigma(1)},\dots,v_{\sigma(k)})
\eta(v_{\sigma(k+1)},\dots,v_{\sigma(k+l)})
\xi(v_{\sigma(k+l+1)},\dots,v_{\sigma(k+l+m)}),
\]
where \(\operatorname{Sh}(k,l,m)\) preserves the order within each of the three blocks. Thus the product is associative.

Interchanging the blocks of lengths \(k\) and \(l\) establishes a bijection between \(\operatorname{Sh}(k,l)\) and \(\operatorname{Sh}(l,k)\). This interchange requires \(kl\) transpositions, so the sign changes by the factor \((-1)^{kl}\). This proves
\[
\omega\wedge\eta=(-1)^{kl}\eta\wedge\omega.
\]

Finally, repeatedly applying the shuffle formula to \(1\)-covectors gives
\[
(\omega^1\wedge\cdots\wedge\omega^k)(v_1,\dots,v_k)
=
\displaystyle\sum_{\sigma\in S_k}\operatorname{sgn}(\sigma)
\prod_{j=1}^k\omega^j(v_{\sigma(j)}),
\]
which is the permutation expansion of the determinant
\(\det(\omega^j(v_i))\).
\end{proof}
 \begin{definition}\label{def:nociones-fundamentales-haz-tensores-alternantes}\index{bundle of alternating tensors@bundle of alternating \(k\)-tensors}
 Let $M$ be a smooth manifold with or without boundary. We define the bundle of alternating
 $k$-tensors by
 \[\Lambda^{k}(T^{*}M):=\coprod_{p\in M}\Lambda^{k}(T_{p}^{*}M).\]
 \end{definition}
 \begin{proposition}\label{prop:nociones-fundamentales-haz-vectorial-suave-rango}
 Let $M$ be a smooth manifold of dimension $n$ with or without boundary.
 Then $\Lambda^{k}(T^{*}M)$ is a smooth vector bundle of rank
 ${n\choose k}$ over $M$.
 \end{proposition}
 \begin{proof}
 Let $(U,x^1,\dots,x^n)$ be a smooth chart on $M$. For each strictly increasing multi-index
 $I=(i_1<\cdots<i_k)$, set
 $\mathbf{d}x^I:=\mathbf{d}x^{i_1}\wedge\cdots\wedge \mathbf{d}x^{i_k}$.
 Proposition~\ref{prop:nociones-fundamentales-espacio-vectorial-dimension-finita-base-coleccion},
 applied at each point $p\in U$, shows that
 \[
 \bigl(\mathbf{d}x^I|_p\bigr)_{|I|=k}
 \]
 is a basis for $\Lambda^k(T_p^*M)$. Thus the map
 \[
 \Lambda^k(T^*M)|_U\longrightarrow U\times\mathbb R^{\binom nk},
 \qquad
 \displaystyle\sum_{|I|=k}a_I\,\mathbf{d}x^I|_p
 \longmapsto
 \bigl(p,(a_I)_{|I|=k}\bigr),
 \]
 is a local trivialization that is linear on the fibers.

 If $(V,y^1,\dots,y^n)$ is another chart, then on $U\cap V$ each
 $\mathbf{d}y^J$ can be expressed as
 \[
 \mathbf{d}y^J
 =
 \displaystyle\sum_{|I|=k}
 \det\!\left(\frac{\partial y^{j_a}}{\partial x^{i_b}}\right)_{a,b=1}^k
 \mathbf{d}x^I.
 \]
 The coefficients are smooth; consequently, the transitions between these trivializations are smooth and linear on the fibers. This gives
 $\Lambda^k(T^*M)$ the structure of a smooth vector bundle, and the number of increasing indices $I$ is $\binom nk$.
 \end{proof}
 \begin{definition}\label{def:nociones-fundamentales-forma-diferencial}\index{differential form}
 Let $M$ be a smooth manifold with or without boundary. A section of
 $\Lambda^{k}(T^{*}M)$ is called a \textbf{differential $k$-form} or simply a \textbf{$k$-form}. The integer $k$ is called the
 \textbf{degree} of the form. We denote the vector space of smooth $k$-forms by \[\Omega^{k}(M):=\Gamma(\Lambda^{k}(T^{*}M)).\]
 \end{definition}
 \begin{remark}\label{obs:nociones-fundamentales-tensor-covariante-alternante}
 Every $1$-covariant tensor is alternating, so $\Omega^{1}(M)=\Gamma(T^{(0,1)}(TM))$.
 \end{remark}
 \begin{proposition}\label{pullback de top degree forms}Let $F\colon M\longrightarrow N$ be a smooth map between $n$-manifolds with or without boundary. If $(x^{i})$ and $(y^{i})$ are smooth coordinates on $U\subseteq M$ and $V\subseteq N$, respectively, and $u\in C^\infty(V)$, then the following holds on $U\cap F^{-1}(V)$:
 \[F^{*}(u\mathbf{d}y^{1}\wedge\cdots\wedge \mathbf{d}y^{n})=(u\circ F)(\det DF)\mathbf{d}x^{1}\wedge \cdots \wedge \mathbf{d}x^{n},\] where $DF$ denotes the Jacobian matrix of $F$ in these coordinates.

 \end{proposition}
\begin{proof}
By the multiplicativity of pullback,
\[
F^*(u\,\mathbf d y^1\wedge\cdots\wedge\mathbf d y^n)
=(u\circ F)\,d(y^1\circ F)\wedge\cdots\wedge d(y^n\circ F).
\]
In the coordinates \(x\),
\[
d(y^a\circ F)=\displaystyle\sum_{b=1}^n
\frac{\partial(y^a\circ F)}{\partial x^b}\,\mathbf d x^b.
\]
Expanding the exterior product, the terms with repeated indices vanish. The sum of the remaining terms, indexed by
\(\sigma\in S_n\), is
\[
\det\left(\frac{\partial(y^a\circ F)}{\partial x^b}\right)_{a,b=1}^n
\mathbf d x^1\wedge\cdots\wedge\mathbf d x^n,
\]
by the Leibniz formula for the determinant. The matrix in parentheses is
\(DF\) in the chosen charts.
\end{proof}
The exterior derivative extends the differential of functions while preserving its intrinsic character.
\begin{theorem}[Existence and uniqueness of the exterior derivative]\label{Existencia y unicidad de la derivada exterior}\index{existence and uniqueness of the exterior derivative}
 Let $M$ be a smooth manifold with or without boundary. For each
 $k\in\mathbb{N}_{0}$, there exists a unique linear map
 \[
 d\colon\Omega^{k}(M)\longrightarrow\Omega^{k+1}(M)
 \]
 that, together with the operators in the other degrees, satisfies the following properties.
 \begin{enumerate}[label=(\alph*)]
 \item If $f\in\Omega^{0}(M)=C^{\infty}(M)$, then $df$ is the differential of $f$, that is, $df(X)=X(f)$.
 \item If $\boldsymbol{\omega}\in\Omega^{k}(M)$ and $\boldsymbol{\eta}\in\Omega^{l}(M)$, then
 \[
 d(\boldsymbol{\omega}\wedge\boldsymbol{\eta})
 =d\boldsymbol{\omega}\wedge\boldsymbol{\eta}+(-1)^{k}\boldsymbol{\omega}\wedge d\boldsymbol{\eta}.
 \]
 \item $d\circ d=0$.
 \item If $F\colon M\longrightarrow N$ is a smooth map and
 $\boldsymbol{\omega}\in\Omega^{k}(N)$, then
 \[
 d(F^{*}\boldsymbol{\omega})=F^{*}(d\boldsymbol{\omega}).
 \]
 \end{enumerate}
 In a smooth chart with coordinates $(x^{1},\dots,x^{n})$, its expression is
 \begin{equation}
 \label{eq:expresion-local-derivada-exterior}
 d\left(\displaystyle\sum_{i_{1}<\cdots<i_{k}}
 \boldsymbol{\omega}_{i_{1}\dots i_{k}}\,
 \mathbf{d}x^{i_{1}}\wedge\cdots\wedge \mathbf{d}x^{i_{k}}\right)
 =
 \displaystyle\sum_{i_{1}<\cdots<i_{k}}
 d\boldsymbol{\omega}_{i_{1}\dots i_{k}}\wedge
 \mathbf{d}x^{i_{1}}\wedge\cdots\wedge \mathbf{d}x^{i_{k}}.
 \end{equation}
\end{theorem}

\begin{proof}
We first construct the operator by an intrinsic formula. For
$\boldsymbol{\omega}\in\Omega^{k}(M)$ and local vector fields
$X_{0},\dots,X_{k}$, define
\begin{equation}
\label{eq:formula-intrinseca-derivada-exterior}
\begin{aligned}
(d\boldsymbol{\omega})(X_{0},\dots,X_{k})
={}&\displaystyle\sum_{i=0}^{k}(-1)^{i}
X_{i}\bigl(\boldsymbol{\omega}(X_{0},\dots,\widehat X_{i},\dots,X_{k})\bigr)\\
&+\displaystyle\sum_{0\leq i<j\leq k}(-1)^{i+j}
\boldsymbol{\omega}\bigl([X_{i},X_{j}],X_{0},\dots,
\widehat X_{i},\dots,\widehat X_{j},\dots,X_{k}\bigr).
\end{aligned}
\end{equation}
In the second sum, the remaining fields appear after the bracket in their original order. For $k=0$, the second sum is empty and the formula reduces to $(df)(X)=X(f)$.

The expression is linear over $\mathbb{R}$ in $\boldsymbol{\omega}$ and in each field. The apparent differential dependence on the $X_i$ cancels, and the result is in fact linear over $C^{\infty}(M)$ in each argument.
Let us verify this for the argument with index $r$. Replacing $X_r$ by
$fX_r$ multiplies every term that does not differentiate $f$ by a factor of $f$.
For $i<r$, the additional term from the first sum is
\[
(-1)^iX_i(f)\,
\boldsymbol{\omega}(X_0,\dots,\widehat X_i,\dots,X_r,\dots,X_k).
\]
In the second sum, the term corresponding to the pair $(i,r)$ contains
\[
[X_i,fX_r]=f[X_i,X_r]+X_i(f)X_r.
\]
Moving $X_r$ from the first argument of $\boldsymbol{\omega}$ to the position it occupied in the preceding expression requires $r-1$ transpositions. Its additional contribution therefore has sign
$(-1)^{i+r}(-1)^{r-1}=-(-1)^i$ and cancels the preceding term. If
$j>r$, the identity
\[
[fX_r,X_j]=f[X_r,X_j]-X_j(f)X_r
\]
yields, after moving $X_r$ through $r$ arguments, the sign
$(-1)^{r+j+1}(-1)^r=-(-1)^j$, which cancels the term in which $X_j$ differentiates
$f$ in the first sum. No other derivatives of $f$ appear. This proves
$C^{\infty}(M)$-linearity in the argument $r$, and $r$ was arbitrary.

Consequently, the value of
\eqref{eq:formula-intrinseca-derivada-exterior} at a point depends only on the values of the fields at that point. In particular, it is independent of the chosen local extensions. Thus the formula defines a well-defined
$(k+1)$-covariant tensor field.

We now show that this tensor is alternating and smooth. In a chart with coordinates $(x^1,\dots,x^n)$, write
\[
\boldsymbol{\omega}=
\displaystyle\sum_{1\leq i_1<\cdots<i_k\leq n}
\omega_{i_1\dots i_k}
\mathbf{d}x^{i_1}\wedge\cdots\wedge \mathbf{d}x^{i_k}.
\]
Since coordinate vector fields commute, evaluating
\eqref{eq:formula-intrinseca-derivada-exterior} on them gives
\[
d\boldsymbol{\omega}=
\displaystyle\sum_{1\leq i_1<\cdots<i_k\leq n}
\displaystyle\sum_{a=1}^{n}
\frac{\partial\omega_{i_1\dots i_k}}{\partial x^a}\,
\mathbf{d}x^a\wedge \mathbf{d}x^{i_1}\wedge\cdots\wedge \mathbf{d}x^{i_k}.
\]
Both sides are tensors and agree on every tuple of coordinate vector fields, so they agree throughout the chart. The right-hand side is a smooth differential form; hence the tensor constructed is alternating and smooth. Moreover, this equality is precisely
\eqref{eq:expresion-local-derivada-exterior}. Since it was derived from the intrinsic formula, the expressions obtained in different charts agree on their overlaps. This proves that the operator is globally well defined and independent of coordinates.

The graded Leibniz rule can be checked in a chart. For $d,q\in\mathbb N_0$, write $\mathcal J_d^q:=\{(i_1,\dots,i_q)\in\{1,\dots,d\}^q\mid i_1<\cdots<i_q\}$, with $\mathcal J_d^0=\{\varnothing\}$. If
$\displaystyle \boldsymbol{\omega}=\displaystyle\sum_{I\in\mathcal J_{n}^k} a_I \mathbf{d}x^I$ has degree $k$ and
$\displaystyle \boldsymbol{\eta}=\displaystyle\sum_{J\in\mathcal J_n^\ell} b_J \mathbf{d}x^J$ has degree $\ell$, the local expression and the product rule for functions give
\[
\begin{aligned}
d(\boldsymbol{\omega}\wedge\boldsymbol{\eta})
&=\displaystyle\sum_{\substack{I\in\mathcal J_n^k\\J\in\mathcal J_n^\ell}}d(a_Ib_J)\wedge \mathbf{d}x^I\wedge \mathbf{d}x^J\\
&=\displaystyle\sum_{\substack{I\in\mathcal J_n^k\\J\in\mathcal J_n^\ell}}da_I\wedge \mathbf{d}x^I\wedge b_J\mathbf{d}x^J
 +\displaystyle\sum_{\substack{I\in\mathcal J_n^k\\J\in\mathcal J_n^\ell}}a_I\,db_J\wedge \mathbf{d}x^I\wedge \mathbf{d}x^J\\
&=d\boldsymbol{\omega}\wedge\boldsymbol{\eta}+(-1)^k\boldsymbol{\omega}\wedge d\boldsymbol{\eta}.
\end{aligned}
\]
The sign in the last equality arises from moving the $1$-form $db_J$ past the $k$ factors of $\mathbf{d}x^I$.

Applying the local expression twice, we obtain
\[
d^{2}\boldsymbol{\omega}
=\displaystyle\sum_{I\in\mathcal J_{n}^k}\displaystyle\sum_{a,b=1}^{n}
\frac{\partial^{2}a_I}{\partial x^{b}\partial x^{a}}
\mathbf{d}x^{b}\wedge \mathbf{d}x^{a}\wedge \mathbf{d}x^{I}.
\]
The terms with $a=b$ vanish. Terms with distinct indices cancel in pairs, since mixed partial derivatives commute and
$\mathbf{d}x^{b}\wedge \mathbf{d}x^{a}=-\mathbf{d}x^{a}\wedge \mathbf{d}x^{b}$. Therefore $d^{2}=0$.

To prove naturality, let $F\colon M\to N$ be smooth and, in a chart on
$N$, let $\displaystyle \boldsymbol{\omega}=\displaystyle\sum_{I\in\mathcal J_{\dim N}^k} a_I\,\mathbf{d}y^I$. The chain rule gives
\[
d(a_I\circ F)=F^*(da_I),
\qquad F^*(\mathbf{d}y^j)=d(y^j\circ F).
\]
Using the graded Leibniz rule and $d^2=0$, we obtain
\[
\begin{aligned}
d(F^*\boldsymbol{\omega})
&=d\left(
\displaystyle\sum_{I\in\mathcal J_{\dim N}^k}(a_I\circ F)\,
d(y^{i_1}\circ F)\wedge\cdots\wedge d(y^{i_k}\circ F)
\right)\\
&=\displaystyle\sum_{I\in\mathcal J_{\dim N}^k} F^*(da_I)\wedge F^*(\mathbf{d}y^I)
=F^*(d\boldsymbol{\omega}).
\end{aligned}
\]
The identity is local around each point and therefore holds on all of $M$.

If $M$ has boundary, the charts used are relatively open in the half-space. Their functions, forms, and fields admit the local extensions specified in the manuscript's definition of smoothness. The intrinsic formula is independent of these extensions, and the coordinate expression above is obtained by restriction. All the arguments, including those using mixed derivatives, remain valid up to the boundary.

Only uniqueness remains. Let $D$ be another family of operators on
$M$, in all degrees, satisfying (a)--(c). First, $D$ is local. If
$\boldsymbol\omega$ vanishes near $p$, choose a smooth function $\chi$
with compact support in that neighborhood and equal to one near $p$.
Since $\chi\boldsymbol\omega=0$, the graded Leibniz rule gives
$0=d\chi\wedge\boldsymbol\omega+\chi D\boldsymbol\omega$; evaluating
at $p$ yields $(D\boldsymbol\omega)_p=0$. Hence forms that agree near
$p$ have the same image under $D$ at $p$.

A form on a coordinate domain can be extended near any of its points
to a global form by multiplication by a cutoff and extension by zero.
Locality makes the resulting local operator independent of this
extension, and (a)--(c) remain valid locally. There,
$D(\mathbf{d}x^i)=D^2x^i=0$, and the graded Leibniz rule implies $D(\mathbf{d}x^I)=0$. Therefore
\[
D\left(\displaystyle\sum_{I\in\mathcal J_{n}^k} a_I \mathbf{d}x^I\right)
=\displaystyle\sum_{I\in\mathcal J_{n}^k} da_I\wedge \mathbf{d}x^I,
\]
which agrees with \eqref{eq:expresion-local-derivada-exterior}. Thus
$D=d$ in every chart and hence on the entire manifold.
\end{proof}
\chapter{Riemannian manifolds}
\label{cap:variedades-riemannianas}

A smooth manifold has local coordinates, but by itself carries no intrinsic notion of length, angle, or volume. Introducing these quantities requires choosing an inner product on each tangent space that depends smoothly on the point. This structure is a Riemannian metric.

The metric defines the length of curves, the Riemannian distance, and the volume measure. Together with the Levi--Civita connection, it also allows us to construct the gradient, divergence, Laplacian, and covariant derivatives that will appear in the analytical chapters. Riemannian geometry therefore provides the intrinsic language in which we will formulate the Sobolev spaces, variational problems, and differential equations studied later.

\begin{figura}[htbp]

    \begin{minipage}[t]{0.30\textwidth}
        \centering
        \includegraphics[height=3.4cm]{500px-Carl_Friedrich_Gauss_1840_by_Jensen.jpg}

        \smallskip
        \footnotesize Carl Friedrich Gauss\\
        (1777--1855)
    \end{minipage}
    \hfill
    \begin{minipage}[t]{0.30\textwidth}
        \centering
        \includegraphics[height=3.4cm]{riemann.jpeg}

        \smallskip
        \footnotesize Bernhard Riemann\\
        (1826--1866)
    \end{minipage}
    \hfill
    \begin{minipage}[t]{0.30\textwidth}
        \centering
        \includegraphics[height=3.4cm]{JH_Poincare.jpg}

        \smallskip
        \footnotesize Henri Poincaré\\
        (1854--1912)
    \end{minipage}

    \caption[Gauss, Riemann, and Poincaré]{Gauss, Riemann, and Poincaré were central figures in the development of modern geometry. Gauss's work on the intrinsic geometry of surfaces paved the way for Riemannian geometry, introduced by Riemann to study spaces of arbitrary dimension. Poincaré, in turn, advanced the global study of these spaces and laid some of the foundations of modern topology.}
    \label{fig:gauss-riemann-poincare}
\end{figura}

\begin{semblanzaHistorica}{From Gauss's surfaces to Riemann's manifolds}
Gauss discovered that the curvature of a surface can be calculated from the metric inherited by its coordinates, without making direct use of the ambient space. Riemann turned this observation into a much broader program: the study of spaces of arbitrary dimension from an infinitesimal rule for measuring lengths. Poincaré later showed that local geometry must coexist with global invariants. This chapter follows that evolution: it begins with the metric, constructs the connection and curvature, and ends with the tools of integration and the operators needed to carry out analysis on the manifold.
\end{semblanzaHistorica}

\section{Riemannian metrics}
A Riemannian metric adds an inner product on each tangent space to the smooth structure, allowing many constructions from Euclidean analysis to be transferred to the geometric setting.

\begin{definition}\label{def:variedades-riemannianas-variedad}\index{Riemannian metric}\index{manifold!Riemannian}\glsadd{variedad-riemanniana}\glsadd{metrica-riemanniana}
 Let $M$ be a smooth manifold with or without boundary.
 A \textit{Riemannian metric on $M$} is a smooth tensor field
 \[
 \mathbf{g} \in \Gamma\bigl(T^{(0,2)}(TM)\bigr)
 \]
 that is symmetric and positive definite at every point $p\in M$.
\end{definition}

\begin{definition}\label{def:variedades-riemannianas-variedad-riemanniana-con-frontera}\index{manifold!Riemannian with boundary}

A \textbf{Riemannian manifold with or without boundary} is a pair $(M,\mathbf{g})$, where $M$ is a smooth manifold with or without boundary and $\mathbf{g}$ is a Riemannian metric on $M$.
\end{definition}
\begin{example}\label{ej:variedades-riemannianas-metrica-euclidiana}
 $\mathbb{R}^{n}$ is a Riemannian manifold without boundary with the \textbf{Euclidean metric} $\overline{\mathbf{g}}$, given in standard coordinates by \[\overline{\mathbf{g}}=\delta_{ij}\,\mathbf{d}x^{i}\otimes\mathbf{d}x^{j}\] or by \[\overline{\mathbf{g}}=(\mathbf{d}x^{1})^{2}+\dots+(\mathbf{d}x^{n})^{2}.\]
 If $v,w\in T_{p}\mathbb{R}^{n}$, $\overline{\mathbf{g}}_{p}(v,w)=\delta_{ij}v^{i}w^{j}=\displaystyle\sum_{i=1}^{n}v^{i}w^{i}=v\cdot w$, which is the usual inner product on $\mathbb{R}^{n}\cong T_{p}\mathbb{R}^{n}$.
\end{example}
\begin{remark}\label{obs:variedades-riemannianas-producto-interno}
 Given a Riemannian manifold $(M,\mathbf{g})$ with or without boundary, the Riemannian metric $\mathbf{g}$ defines an \textbf{inner product} on each tangent space $T_{p}M$ for each $p\in M$. For $\mathbf{X},\mathbf{Y}\in \mathfrak{X}(M)$, we sometimes denote $\mathbf{g}(\mathbf{X},\mathbf{Y})\in C^{\infty}(M)$ by $\mathbf{g}(\mathbf{X},\mathbf{Y}):=\langle \mathbf{X},\mathbf{Y} \rangle$. Pointwise, we have $\langle \mathbf{X}_{p},\mathbf{Y}_{p}\rangle=\mathbf{g}_{p}(\mathbf{X}_{p
 },\mathbf{Y}_{p})\in\mathbb{R}$. We can also define a \textbf{norm} on $T_{p}M$ by $\|v\|_{\mathbf{g}}=\sqrt{\langle v,v\rangle_{\mathbf{g}}}$.

 Finally, the concept of angle between two vectors extends from $\mathbb{R}^{n}$ to $M$ as follows:

 The angle between two vectors $v,w\in T_{p}M\setminus \{0\}$ is the unique $\theta\in [0,\pi]$ satisfying \[\cos(\theta)=\frac{\langle v,w\rangle_{\mathbf{g}}}{\|v\|_{\mathbf{g}}\|w\|_{\mathbf{g}}}.\] Two vectors $v,w\in T_{p}M$ are said to be \textbf{orthogonal} if $\langle v,w \rangle_{\mathbf{g}}=0$.
\end{remark}
\begin{definition}[Pointwise order on symmetric forms and comparison of metrics]
\label{def:orden-formas-metricas}
Let $\mathbf a$ and $\mathbf b$ be two symmetric tensor fields of type
$(0,2)$, and let $A\subseteq M$. The notation
$\mathbf a\leq\mathbf b$ on $A$ means that
\[
 \mathbf a_p(v,v)\leq\mathbf b_p(v,v)
 \qquad\text{for every }p\in A\text{ and every }v\in T_pM.
\]
Equivalently, $(\mathbf b-\mathbf a)_p(v,v)\geq0$ for the same
$p$ and $v$; that is, $\mathbf b-\mathbf a$ is positive semidefinite on each fiber over $A$. This is the meaning of an inequality between symmetric tensors of type $(0,2)$ in this book.

In particular, if $\mathbf g$ and $\widetilde{\mathbf g}$ are Riemannian metrics and $0<c\leq C<\infty$, writing
\[
 c\mathbf g\leq\widetilde{\mathbf g}\leq C\mathbf g
 \qquad\text{on }A
\]
is shorthand for the numerical inequalities
\[
 c\,\mathbf g_p(v,v)\leq\widetilde{\mathbf g}_p(v,v)
 \leq C\,\mathbf g_p(v,v)
 \qquad(p\in A,\ v\in T_pM).
\]
Since $\mathbf g_p(v,v)=|v|_{\mathbf g}^2$, this is exactly equivalent to
\begin{equation}
\label{eq:orden-metricas-normas-vectores}
 \sqrt c\,|v|_{\mathbf g}\leq|v|_{\widetilde{\mathbf g}}
 \leq\sqrt C\,|v|_{\mathbf g}
 \qquad(p\in A,\ v\in T_pM).
\end{equation}
The constants are common to all points of $A$ and all vectors in their fibers. When such constants exist, we say that the metrics are uniformly equivalent on $A$.
\end{definition}

\begin{remark}
\label{obs:interpretacion-orden-metricas}
The comparison is made by evaluating both metrics on the same vector, repeated in both arguments. If $v=\displaystyle\sum_{i=1}^{n} v^i\partial_i$, its expression in a chart is
\[
 c\displaystyle\sum_{i,j=1}^{n}g_{ij}(p)v^iv^j
 \leq\displaystyle\sum_{i,j=1}^{n}\widetilde g_{ij}(p)v^iv^j
 \leq C\displaystyle\sum_{i,j=1}^{n}g_{ij}(p)v^iv^j
 \qquad\text{for every }(v^1,\ldots,v^n)\in\mathbb R^n.
\]
The matrix entries, considered individually, need not satisfy these inequalities. Nor is an inequality asserted for $\mathbf g_p(v,w)$ and
$\widetilde{\mathbf g}_p(v,w)$ with arbitrary $v$ and $w$.
Diagonal evaluation is independent of the chart and determines the symmetric form through the polarization identity
\[
 \mathbf g_p(v,w)
 =\frac14\bigl(\mathbf g_p(v+w,v+w)-\mathbf g_p(v-w,v-w)\bigr).
\]

Comparison constants exist on each finite-dimensional fiber; uniformity in the base point is an additional assertion.
On a compact set, it follows from continuity of the metrics, as will be proved in Lemma~\ref{lema cg} and, for arbitrary vector bundles, in Lemma~\ref{lema equivalencia metricas fibradas}.
The same diagonal evaluation convention is used to compare bundle metrics, with $v\in\mathbf E_p$; in the complex case, we compare the real values of the Hermitian forms on $(v,v)$.
\end{remark}

We can use the pullback of tensor fields, defined in Chapter 1, to obtain a Riemannian metric from smooth maps.
\begin{proposition}\label{pullback de métricas riemannianas}
 Let $M$ and $N$ be smooth manifolds with or without boundary, let
 $F\colon M\longrightarrow N$ be a smooth map, and let $\mathbf{g}$ be a Riemannian metric on $N$. Then $F^{*}\mathbf{g}$ is a Riemannian metric on $M$ if and only if $F$ is an immersion.
\end{proposition}
\begin{proof}
 By Proposition~\ref{propiedades del pullback}, we know that $F^{*}\mathbf{g}\in \Gamma(T^{(0,2)}(TM))$, that is, $F^{*}\mathbf{g}$ is a smooth $2$-covariant tensor field. The definition of pullback shows that this tensor is symmetric; it only remains to determine when it is positive definite at every point.

First suppose that $F^{*}\mathbf{g}$ is a Riemannian metric on $M$. We show that $F$ is an immersion. Take $p\in M$ and let $v\in T_{p}M$ with $v\neq 0$. By Definition~\ref{def:inmersion-submersion}, we must show that $dF_{p}$ is injective. If $dF_{p}(v)=0$, then \[(F^{*}\mathbf{g})_{p}(v,v)=g_{F(p)}(dF_{p}(v),dF_{p}(v))=0\] contradicting the positive definiteness of $F^{*}\mathbf{g}$ on $p$. It follows that $dF_p(v)\neq0$ for every $v\neq0$; thus $dF_p$ is injective and $F$ is an immersion.

Conversely, suppose that $F$ is an immersion. If $p\in M$ and $v\in T_pM$ satisfy
 \[
0=(F^{*}\mathbf{g})_{p}(v,v)=\mathbf{g}_{F(p)}(dF_{p}(v),dF_{p}(v)),
 \]
then, since $\mathbf{g}_{F(p)}$ is positive definite, the preceding equality implies $dF_p(v)=0$. The injectivity of $dF_p$ then gives $v=0$, so $F^*\mathbf{g}$ is positive definite.
\end{proof}
Every smooth manifold, with or without boundary, admits at least one Riemannian metric.
\begin{proposition}[Existence of Riemannian metrics]\label{prop:variedades-riemannianas-existencia-de-metricas-riemannianas}\index{existence of Riemannian metrics} Every smooth manifold with or without boundary admits a Riemannian metric.

\end{proposition}
\begin{proof}
 Let $M$ be a smooth manifold with or without boundary. Let $(U_{\alpha},\phi_{\alpha})_{\alpha\in J}$ be a family of charts such that $(U_{\alpha})_{\alpha\in J}$ is an open cover of $M$. By Theorem~\ref{particionesdelaunidad}, we can find a partition of unity $(\psi_{\alpha})_{\alpha\in J}$ subordinate to the cover $(U_{\alpha})_{\alpha\in J}$. By Proposition~\ref{pullback de métricas riemannianas}, in each coordinate chart, $\mathbf{g}_{\alpha}=\phi_{\alpha}^{*}\overline{\mathbf{g}}$ defines a Riemannian metric on $U_{\alpha}$ whose coordinate expression is $\delta_{ij}\,\mathbf{d}x^{i}\otimes \mathbf{d}x^{j}$. Since $\operatorname{supp}(\psi_\alpha)\subseteq U_\alpha$, the tensor $\psi_\alpha \mathbf{g}_\alpha$ extends smoothly by zero outside $U_\alpha$. Define
 \[
 \mathbf{g}:=\displaystyle\sum_{\alpha\in J}\psi_{\alpha}\mathbf{g}_{\alpha},
 \]
 where each summand is understood to have this extension. The family of supports is locally finite, so the sum is finite around each point and defines a smooth tensor field.

The field $\mathbf{g}$ is symmetric because it is a finite sum of symmetric tensors in a neighborhood of each point.

 To see that it is positive definite, let $p\in M$ and $v\in T_pM\setminus\{0\}$. For each index with $\psi_\alpha(p)>0$, we have $p\in U_\alpha$ and $(\mathbf{g}_\alpha)_p(v,v)>0$, while all the coefficients $\psi_\alpha(p)$ are nonnegative. Since their sum is one, there exists $\beta$ such that $\psi_\beta(p)>0$. Therefore,
 \[
 \mathbf{g}_p(v,v)=\displaystyle\sum_{\alpha\in J}\psi_\alpha(p)(\mathbf{g}_\alpha)_p(v,v)>0.
 \]
\end{proof}
If $(\mathbf{E}_{1},\dots,\mathbf{E}_{n})$ is a smooth local frame for $TM$ defined on an open subset $U\subseteq M$ and $(\boldsymbol{\varepsilon}^{1},\dots,\boldsymbol{\varepsilon}^{n})$ is its dual coframe, we can write $\mathbf{g}$ locally on $U$ as
\[
\mathbf{g}=g_{ij}\,\boldsymbol{\varepsilon}^{i}\otimes\boldsymbol{\varepsilon}^{j},
\qquad
g_{ij}(p):=\langle \mathbf{E}_i|_p,\mathbf{E}_j|_p\rangle_{\mathbf{g}}.
\]
These functions are called the \textit{\textbf{coefficients of the metric $\mathbf{g}$ with respect to the local frame $(\mathbf{E}_{1},\dots,\mathbf{E}_{n})$}}. We say that $(\mathbf{E}_{1},\dots,\mathbf{E}_{n})$ is an \textbf{\textit{orthonormal local frame}} if the vectors $\mathbf{E}_{1}|_{p},\dots,\mathbf{E}_{n}|_{p}$ form an orthonormal basis of $T_{p}M$ for every $p\in U$. Equivalently, a local frame is orthonormal if and only if $\langle \mathbf{E}_{i},\mathbf{E}_{j}\rangle_{\mathbf{g}}=\delta_{ij}$.

 \begin{proposition}[Existence of orthonormal local frames]\label{marco ortonrmal existencia}\index{existence of orthonormal local frames} Let $(M,\mathbf{g})$ be a Riemannian manifold with or without boundary. If $(\mathbf{X}_{1},\dots,\mathbf{X}_{n})$ is a smooth local frame for $TM$ defined on $U\subseteq M$, then there exists a smooth orthonormal local frame $(\mathbf{E}_{1},\dots,\mathbf{E}_{n})$ defined on $U$ such that
 \[
 \operatorname{span}\{\mathbf{E}_1|_p,\dots,\mathbf{E}_k|_p\}
 =\operatorname{span}\{\mathbf{X}_1|_p,\dots,\mathbf{X}_k|_p\}
 \]
 for each $k\in\{1,\dots,n\}$ and each $p\in U$. In particular, for each $p\in M$, there exists an orthonormal local frame defined on some neighborhood of $p$.

\end{proposition}
\begin{proof}
 By the Gram--Schmidt algorithm, for each $p\in U$, the vectors
 $\mathbf{E}_1|_p,\dots,\mathbf{E}_n|_p$ can be chosen recursively so that
 \[
 \operatorname{span}\{\mathbf{E}_1|_p,\dots,\mathbf{E}_k|_p\}
 =\operatorname{span}\{\mathbf{X}_1|_p,\dots,\mathbf{X}_k|_p\}
 \]
 for every $k$. Explicitly,
 $\mathbf{E}_{1}=\frac{\mathbf{X}_{1}}{\|\mathbf{X}_{1}\|_{\mathbf{g}}}$ and
 \[
 \mathbf{E}_{j}=\displaystyle\frac{\mathbf{X}_{j}-\displaystyle\sum_{i=1}^{j-1}\langle \mathbf{X}_{j},\mathbf{E}_{i}\rangle_{\mathbf{g}} \mathbf{E}_{i}}
 {\left\|\mathbf{X}_{j}-\displaystyle\sum_{i=1}^{j-1}\langle \mathbf{X}_{j},\mathbf{E}_{i}\rangle_{\mathbf{g}} \mathbf{E}_{i}\right\|_{\mathbf{g}}},
 \qquad 2\leq j\leq n.
 \]
 Each denominator is a smooth, strictly positive function because
 $(\mathbf{X}_1,\dots,\mathbf{X}_n)$ is a frame. Therefore the fields $\mathbf{E}_i$ are smooth;
 the Gram--Schmidt identities show that they form an orthonormal frame and satisfy the stated equalities of subspaces.
\end{proof}

Riemannian metrics allow us to introduce a bundle isomorphism between $TM$ and $T^{*}M$ as follows:

\begin{definition}[Musical isomorphisms]\label{isomorfismos musicales}\index{musical isomorphisms}
Let $(M,\mathbf{g})$ be a Riemannian manifold with or without boundary. For each point $p\in M$, define the maps
\[
\mathbf{g}^\flat_p\colon T_pM \longrightarrow T_p^*M,
\qquad
X_p \mapsto \mathbf{g}_p(X_p,\cdot),
\]
and
\[
\mathbf{g}^\sharp_p\colon T_p^*M \longrightarrow T_pM,
\qquad
\omega_p \mapsto \text{the unique $X_p$ such that } \mathbf{g}_p(X_p,\cdot)=\omega_p.
\]
These maps are inverses of one another and define bundle isomorphisms
\[
\mathbf{g}^\flat\colon TM \longrightarrow T^*M,
\qquad
\mathbf{g}^\sharp\colon T^*M \longrightarrow TM,
\]
called the \emph{musical isomorphisms}.
\end{definition}

\begin{remark}\label{obs:variedades-riemannianas-coordenadas-locales-campo-vectorial-forma-inversa}
In local coordinates $(x^1,\dots,x^n)$, if
$\mathbf{X}=X^i\boldsymbol{\partial}_i$ is a vector field, then
\[
\mathbf{X}^\flat = g_{ij}X^i\,\mathbf{d}x^j,
\]
and if $\boldsymbol{\omega}=\omega_i\,\mathbf{d}x^i$ is a $1$-form,
\[
\boldsymbol{\omega}^\sharp = g^{ij}\omega_i\boldsymbol{\partial}_j,
\]
where $(g^{ij})$ is the inverse of the metric matrix $(g_{ij})$.
Thus $^\flat$ is interpreted as “lowering indices” and $^\sharp$ as “raising indices.”
\end{remark}

 \section{Traces of tensors}
Contraction removes one covariant and one contravariant index by means of a sum independent of the basis. This operation will appear repeatedly in the definitions of divergence, Laplacians, curvature, and energy.

 \noindent Given a finite-dimensional vector space $V$ and an endomorphism $T\colon V\longrightarrow V$, if $\beta=\{e_{1},\dots,e_{n}\}$ is a basis for $V$ and $(A_{ij})$ is the matrix of $T$ with respect to the basis $\beta$, we can define $\operatorname{tr}(T):=\displaystyle\sum_{i=1}^{n}A_{i,i}$. This definition is independent of the chosen basis, and we can generalize the construction to mixed tensors:

 \begin{definition}\label{def:variedades-riemannianas-operador-traza}\index{trace operator}
 For the $(1,1)$-tensor $F(\omega^{1},\dots,\omega^{k},\cdot,v_{1},\dots,v_{l},\cdot)\in T^{(1,1)}(V)$, we can calculate its trace (which is independent of the chosen basis) $\operatorname{tr}(F)(\omega^{1},\dots,\omega^{k},v_{1},\dots,v_{l})$. In coordinates, the components of $\operatorname{tr}(F)$ are \[(\operatorname{tr}(F))_{j_{1}\dots j_{l}}^{i_{1}\dots i_{k}}=F^{i_{1}\dots i_{k}m}_{j_{1}\dots j_{l}m}.\] This construction defines an operator, independent of the chosen basis, that we call the \emph{trace operator} $\operatorname{tr}\colon T^{(k+1,l+1)}(V)\longrightarrow T^{(k,l)}(V)$.
 \end{definition}
 The trace operator is also known as \textit{\textbf{contraction}}.

In constructing the trace operator, we set the last upper and lower indices equal, which, in Einstein notation, means summing all these terms. A similar construction can be performed for any pair of indices, provided that one is covariant and the other contravariant. In general, we do not distinguish this operation notationally from the operator $\operatorname{tr}$, so we use only this notation and specify the contraction explicitly when necessary.

 The operator $\operatorname{tr}$ can be applied to any pair consisting of one contravariant and one covariant index. It is a well-defined linear map from
 $T^{(k+1,l+1)}(V)$ to $T^{(k,l)}(V)$. If $\mathbf{F}$ is a tensor field of type
 $(k+1,l+1)$ on a manifold $M$, its contraction is defined pointwise by
 \[
 \operatorname{tr}(\mathbf{F})|_p:=\operatorname{tr}(\mathbf{F}|_p),
 \qquad p\in M.
 \]

 Now let $(M,\mathbf{g})$ be a Riemannian manifold with or without boundary, and let
 $\mathbf{F}$ be a $k$-covariant tensor field. The metric allows us to contract two distinct covariant indices: when positions $r$ and $s$ are contracted, the coordinate expression is
 \[
 g^{i_ri_s}F_{i_1\dots i_r\dots i_s\dots i_k},
 \qquad r\ne s.
 \]
 This expression transforms tensorially and therefore defines a
 $(k-2)$-covariant field independent of coordinates. We denote it by
 $\operatorname{tr}_{\mathbf{g}}(\mathbf{F})$; when more than one pair is possible, the two indices being contracted will be specified explicitly.
\section{Fiberwise inner products of tensor fields}
Given a Riemannian manifold $(M,\mathbf{g})$ with or without boundary, we can define several useful concepts from $\mathbf{g}$. One is the inner product of tensor fields. This is a special case of a more general construction that does not necessarily require a Riemannian structure, known as a fiberwise inner product or bundle metric on a vector bundle $\mathbf{E}$.
\begin{definition}\label{def producto fibrado}\index{bundle metric}\index{vector bundle!bundle metric}\glsadd{metrica-fibrada}
 Let $M$ be a smooth manifold with or without boundary. If
 $\pi\colon \mathbf{E}\longrightarrow M$ is a smooth vector bundle, a
 \textbf{smooth bundle metric on $\mathbf{E}$} is an inner product on each fiber $\mathbf{E}_{p}$ that varies smoothly in the sense that, for any smooth local sections
 $\boldsymbol{\sigma},\boldsymbol{\tau}$ of $\mathbf{E}$,
 $\langle\boldsymbol{\sigma},\boldsymbol{\tau}\rangle$ is a smooth function.
\end{definition}
The following example constructs a bundle metric on a nontrivial vector bundle.
\begin{example}[Bundle metric on the Möbius bundle]\label{haz:mobius-metrica}

Return to the bundle $\mathbf{E}\to\mathbb{S}^{1}$ of
Example~\ref{haz:mobius}. For $p\in\mathbb{S}^{1}$, define an inner product on the fiber $\pi^{-1}(p)$ as follows. If $x_{0}\in\mathbb{R}$ satisfies $\varepsilon(x_{0})=p$, then
\[
\pi^{-1}(p)=\{[(x_{0},t)]\mid t\in\mathbb{R}\},
\]
and we consider the map
\[
j_{p}\colon \pi^{-1}(p)\longrightarrow\mathbb{R},\qquad j_{p}([(x_{0},t)])=t.
\]
This map is well defined since, if $(x_{0},t)\sim(x_{0},t')$, then $t=t'$. We then define
\[
\mathbf{h}_{p}(u,v)=j_{p}(u)j_{p}(v),\qquad u,v\in\pi^{-1}(p).
\]
In this way, each fiber is equipped with an inner product. It remains to show that $\mathbf{h}$ is a smooth bundle metric in the sense of Definition~\ref{def producto fibrado}.

Consider the cover $U_{\pm}=\mathbb{S}^{1}\setminus\{(\pm 1,0)\}$ with trivializations
\[
\Phi_{\pm}\colon \pi^{-1}(U_{\pm})\longrightarrow U_{\pm}\times\mathbb{R},\qquad [(x,y)]\longmapsto(\varepsilon(x),y).
\]
On the overlap $U_{+}\cap U_{-}$, we have
\[
(\Phi_{+}\circ\Phi_{-}^{-1})(p,v)=(p,\tau_{+-}(p)v),\qquad \tau_{+-}(p)\in\{\pm 1\}.
\]
This means that if a local section $\boldsymbol{\sigma}$ is written on $U_{-}$ as $\Phi_{-}(\boldsymbol{\sigma}(p))=(p,s_{-}(p))$ with $s_{-}\colon U_{-}\longrightarrow\mathbb{R}$ smooth, then, for every $p\in U_{+}\cap U_{-}$, the following holds:
\[
\Phi_{+}(\boldsymbol{\sigma}(p))
=\big(\Phi_{+}\circ\Phi_{-}^{-1}\big)\big(\Phi_{-}(\boldsymbol{\sigma}(p))\big)
=\big(\Phi_{+}\circ\Phi_{-}^{-1}\big)\big(p,s_{-}(p)\big)
=(p,\tau_{+-}(p)s_{-}(p)).
\]
Since $\Phi_{+}(\boldsymbol{\sigma}(p))=(p,s_{+}(p))$ by definition, comparing second coordinates yields
\[
s_{+}(p)=\tau_{+-}(p)s_{-}(p).
\]
Similarly, if $\boldsymbol{\tau}$ is another section, then on the overlap
\[
t_{+}(p)=\tau_{+-}(p)t_{-}(p),
\]
where $s_{\pm},t_{\pm}$ are the functions representing $\boldsymbol{\sigma},\boldsymbol{\tau}$ in the trivializations $U_{\pm}$. Consequently,
\[
s_{+}(p)t_{+}(p)=\big(\tau_{+-}(p)s_{-}(p)\big)\big(\tau_{+-}(p)t_{-}(p)\big)=s_{-}(p)t_{-}(p),
\]
and the local expressions for the inner product agree on the intersection.

Finally, if $\boldsymbol{\sigma},\boldsymbol{\tau}$ are smooth local sections over $U_{\pm}$, we have
\[
\Phi_{\pm}(\boldsymbol{\sigma}(p))=(p,s_{\pm}(p)),\qquad \Phi_{\pm}(\boldsymbol{\tau}(p))=(p,t_{\pm}(p)),
\]
and the definition of $\mathbf{h}$ gives
\[
\mathbf{h}_{p}(\boldsymbol{\sigma}(p),\boldsymbol{\tau}(p))=s_{\pm}(p)t_{\pm}(p),
\]
which is a smooth function on $p$.

We conclude that $\mathbf{h}$ defines a smooth bundle metric on $(\mathbf{E},\pi,\mathbb{S}^{1})$, independent of the chosen trivialization and obtained naturally by identifying each fiber with $\mathbb{R}$ through $j_{p}$.
\end{example}

\begin{proposition}\label{metrica en tensores}
 Let $(M,\mathbf{g})$ be a Riemannian manifold with or without boundary. Then there exists a unique smooth bundle metric on each tensor bundle $T^{(k,l)}(TM)$ with the property that, if $\boldsymbol{\alpha}_{1},\dots,\boldsymbol{\alpha}_{k+l},\boldsymbol{\beta}_{1},\dots,\boldsymbol{\beta}_{k+l}$ are vector or covector fields as appropriate, then \[\langle \boldsymbol{\alpha}_{1}\otimes\dots\otimes \boldsymbol{\alpha}_{k+l},\boldsymbol{\beta}_{1}\otimes\dots\otimes\boldsymbol{\beta}_{k+l}\rangle_{\mathbf{g}}=\langle\boldsymbol{\alpha}_{1},\boldsymbol{\beta}_{1}\rangle_{\mathbf{g}}\dots\langle \boldsymbol{\alpha}_{k+l},\boldsymbol{\beta}_{k+l}\rangle_{\mathbf{g}}.\] With this inner product, if $(\mathbf{E}_{1},\dots,\mathbf{E}_{n})$ is an orthonormal local frame for $TM$ and $(\boldsymbol{\varepsilon}^{1},\dots,\boldsymbol{\varepsilon}^{n})$ is the corresponding dual coframe, then the collection of tensor fields \[(\mathbf{E}_{i_{1}}\otimes\dots\otimes \mathbf{E}_{i_{k}}\otimes \boldsymbol{\varepsilon}^{j_{1}}\otimes \dots\otimes \boldsymbol{\varepsilon}^{j_{l}})_{1\leq i_{1},\dots,i_{k},j_{1},\dots,j_{l}\leq n}\] forms an orthonormal local frame on $T^{(k,l)}(TM)$ (orthonormal with respect to this bundle metric). In terms of any local frame, not necessarily orthonormal, this bundle metric satisfies the following: for any $\mathbf{X},\mathbf{Y}\in \mathfrak{X}(M)$, $\boldsymbol{\omega},\boldsymbol{\eta}\in \Omega^{1}(M)$ \[\langle \mathbf{X},\mathbf{Y}\rangle_{\mathbf{g}}=g_{ir}X^{i}Y^{r}, \text{ }\langle \boldsymbol{\omega},\boldsymbol{\eta}\rangle_{\mathbf{g}}=g^{js}\omega_{j}\eta_{s}.\] Moreover, for $\mathbf{F},\mathbf{G}\in\Gamma(T^{(k,l)}(TM))$, \[\langle \mathbf{F},\mathbf{G}\rangle_{\mathbf{g}} =g_{i_{1}r_{1}}\dots g_{i_{k}r_{k}}g^{j_{1}s_{1}}\dots g^{j_{l}s_{l}}F^{i_{1}\dots i_{k}}_{j_{1}\dots j_{l}}G^{r_{1}\dots r_{k}}_{s_{1}\dots s_{l}}.\]

 We sometimes write $\mathbf{g}(\mathbf{F},\mathbf{G}):=\langle \mathbf{F},\mathbf{G}\rangle_{\mathbf{g}}$.
\end{proposition}
\begin{proof}
 Consider $p\in M$, and let $(\mathbf{E}_{1},\dots,\mathbf{E}_{n})$ be a smooth local frame for $TM$ defined on an open set $U\subseteq M$ such that $p\in U$. We can define the bundle metric by \[\langle \mathbf{F},\mathbf{G}\rangle =g_{i_{1}r_{1}}\dots g_{i_{k}r_{k}}g^{j_{1}s_{1}}\dots g^{j_{l}s_{l}}F^{i_{1}\dots i_{k}}_{j_{1}\dots j_{l}}G^{r_{1}\dots r_{k}}_{s_{1}\dots s_{l}}\] on $U$ for any tensor fields $\mathbf{F},\mathbf{G}\in \Gamma(T^{(k,l)}(TM))$. This definition is local. To show that it is well defined, we must first verify that it is independent of the chosen local frame; this will also prove uniqueness:

 Suppose that $(\widetilde{\mathbf{E}}_{1},\dots,\widetilde{\mathbf{E}}_{n})$ is another local frame defined on an open set $V$ such that $p\in V$. On $U\cap V$, the local frames are related by \[\mathbf{E}_{i}=A_{i}^{j}\widetilde{\mathbf{E}}_{j}\] where $(A_{i}^{j})$ is a matrix of smooth functions defined on $U\cap V$. Let $(\boldsymbol{\varepsilon}^{1},\dots,\boldsymbol{\varepsilon}^{n})$ and $(\widetilde{\boldsymbol{\varepsilon}}^{1},\dots,\widetilde{\boldsymbol{\varepsilon}}^{n})$ be the local coframes associated with these frames. Let $(A^{-1})_{i}^{j}$ be the inverse matrix of $(A_{j}^{i})$, where the lower index denotes rows and the upper index columns. We have $\widetilde{\mathbf{E}}_{i}=(A^{-1})_{i}^{j}\mathbf{E}_{j}$, and we can also write $\boldsymbol{\varepsilon}^{i}=B_{m}^{i}\widetilde{\boldsymbol{\varepsilon}}^{m}$, so \[(A^{-1})_{j}^{i}=(A^{-1})_{j}^{k}\delta_{k}^{i}=(A^{-1})_{j}^{k}\boldsymbol{\varepsilon}^{i}(\mathbf{E}_{k})=\boldsymbol{\varepsilon}^{i}((A^{-1})_{j}^{k}\mathbf{E}_{k})=\boldsymbol{\varepsilon}^{i}(\widetilde{\mathbf{E}}_{j})=B_{m}^{i}\widetilde{\boldsymbol{\varepsilon}}^{m}(\widetilde{\mathbf{E}}_{j})=B_{m}^{i}\delta^{m}_{j}=B^{i}_{j}.\] We conclude that $\boldsymbol{\varepsilon}^{i}=(A^{-1})_{j}^{i}\widetilde{\boldsymbol{\varepsilon}}^{j}$. In turn, since $\mathbf{E}_{i}=A_{i}^{j}\widetilde{\mathbf{E}}_{j}$, we deduce that $\widetilde{\boldsymbol{\varepsilon}}^{i}=A_{j}^{i}\boldsymbol{\varepsilon}^{j}$.

 We now verify that the expression is independent of the local frame. Let $\mathbf{F},\mathbf{G}\in \Gamma(T^{(k,l)}(TM))$. On $U\cap V$, we have:
 \[
 g_{i_{1}r_{1}}\dots g_{i_{k}r_{k}}g^{j_{1}s_{1}}\dots g^{j_{l}s_{l}}F^{i_{1}\dots i_{k}}_{j_{1}\dots j_{l}}G^{r_{1}\dots r_{k}}_{s_{1}\dots s_{l}}\] \[=(A_{i_{1}}^{a_{1}}A_{r_{1}}^{b_{1}}\widetilde g_{a_{1}b_{1}})\cdots (A_{i_{k}}^{a_{k}}A_{r_{k}}^{b_{k}}\widetilde g_{a_{k}b_{k}})((A^{-1})^{j_{1}}_{e_{1}}(A^{-1})^{s_{1}}_{f_{1}}\widetilde{g}^{e_{1}f_{1}})\cdots ((A^{-1})^{j_{l}}_{e_{l}}(A^{-1})^{s_{l}}_{f_{l}}\widetilde{g}^{e_{l}f_{l}})\] \[((A^{-1})^{i_{1}}_{c_{1}}\cdots (A^{-1})^{i_{k}}_{c_{k}}A_{j_{1}}^{m_{1}}\cdots A^{m_{l}}_{j_{l}}\widetilde{F}^{c_{1}\dots c_{k}}_{m_{1}\dots m_{l}})((A^{-1})^{r_{1}}_{d_{1}}\cdots (A^{-1})^{r_{k}}_{d_{k}}A_{s_{1}}^{n_{1}}\cdots A^{n_{l}}_{s_{l}}\widetilde{G}^{d_{1}\dots d_{k}}_{n_{1}\dots n_{l}}) \] \[=A_{i_{1}}^{a_{1}}(A^{-1})^{i_{1}}_{c_{1}}\cdots A_{i_{k}}^{a_{k}}(A^{-1})^{i_{k}}_{c_{k}}A_{r_{1}}^{b_{1}}(A^{-1})^{r_{1}}_{d_{1}}\cdots A_{r_{k}}^{b_{k}}(A^{-1})^{r_{k}}_{d_{k}}\widetilde g_{a_{1}b_{1}}\cdots\widetilde g_{a_{k}b_{k}}\]
\[(A^{-1})_{e_{1}}^{j_{1}}A_{j_{1}}^{m_{1}}\cdots (A^{-1})_{e_{l}}^{j_{l}}A_{j_{l}}^{m_{l}}(A^{-1})_{f_{1}}^{s_{1}}A_{s_{1}}^{n_{1}}\cdots(A^{-1})_{f_{l}}^{s_{l}}A_{s_{l}}^{n_{l}}\widetilde{g}^{e_{1}f_{1}}\dots \widetilde{g}^{e_{l}f_{l}}\widetilde{F}^{c_{1}\dots c_{k}}_{m_{1}\dots m_{l}}\widetilde{G}^{d_{1}\dots d_{k}}_{n_{1}\dots n_{l}} \]
 \[=\delta^{a_{1}}_{c_{1}}\cdots \delta^{a_{k}}_{c_{k}}\delta^{b_{1}}_{d_{1}}\cdots \delta^{b_{k}}_{d_{k}}\widetilde g_{a_{1}b_{1}}\cdots\widetilde g_{a_{k}b_{k}}\delta^{m_{1}}_{e_{1}}\cdots\delta^{m_{l}}_{e_{l}}\delta^{n_{1}}_{f_{1}}\cdots \delta^{n_{l}}_{f_{l}}\widetilde{g}^{e_{1}f_{1}}\dots \widetilde{g}^{e_{l}f_{l}}\widetilde{F}^{c_{1}\dots c_{k}}_{m_{1}\dots m_{l}}\widetilde{G}^{d_{1}\dots d_{k}}_{n_{1}\dots n_{l}}\] \[=\widetilde g_{a_{1}b_{1}}\cdots \widetilde g_{a_{k}b_{k}}\widetilde{g}^{e_{1}f_{1}}\cdots \widetilde{g}^{e_{l}f_{l}}\widetilde{F}^{a_{1}\dots a_{k}}_{e_{1}\dots e_{l}}\widetilde{G}^{b_{1}\dots b_{k}}_{f_{1}\dots f_{l}}.\]

 Here the penultimate equality follows from the definitions of matrix multiplication and the inverse matrix, that is, $(A^{-1})^{i}_{k}A^{k}_{j}=\delta^{i}_{j}=A_{j}^{m}(A^{-1})_{m}^{i}$.

 Therefore $\langle \mathbf{F},\mathbf{G}\rangle_{\mathbf{g}}$ is independent of the chosen local frame.

 Now that the inner product is well defined, we see that it is bilinear and symmetric, so it suffices to prove positive definiteness. If $k=l=0$, it is the usual product on $\mathbb R$. If only one of the two indices is zero, the calculation below remains valid after omitting the absent index families and beginning the invertibility argument with the first family present. We may therefore assume that $k,l\geq1$. Let $(\mathbf{E}_{1},\dots,\mathbf{E}_{n})$ be a smooth local frame for $TM$ defined on an open set $U\subseteq M$, and consider $\mathbf{F}\in\Gamma(T^{(k,l)}(TM))$. Since the matrix $(g_{ij})$ is positive definite, the spectral theorem for symmetric matrices gives a positive definite symmetric matrix $R$, hence invertible, such that $(g_{ij})=R^{2}$.

 We denote the entries of $R$ by $R_{ij}$ and those of $R^{-1}$ by $R^{ij}$. In particular, we have $g_{ij}=\displaystyle\sum_{m=1}^{n}R_{im}R_{mj}$ and $g^{ij}=\displaystyle\sum_{m=1}^{n}R^{im}R^{mj}$.

 In the following calculation, we write the sums in Einstein notation explicitly:

 \[\langle \mathbf{F},\mathbf{F}\rangle_{\mathbf{g}}=\displaystyle\sum _{\substack{1\leq i_{1},\dots,i_{k}\leq n\\1\leq j_{1},\dots,j_{l}\leq n\\ 1\leq r_{1},\dots,r_{k}\leq n\\ 1\leq s_{1},\dots,s_{l}\leq n}}g_{i_{1}r_{1}}\cdots g_{i_{k}r_{k}}g^{j_{1}s_{1}}\cdots g^{j_{l}s_{l}}F_{j_{1}\dots j_{l}}^{i_{1}\dots i_{k}}F_{s_{1}\dots s_{l}}^{r_{1}\dots r_{k}}=\] {\tiny \[ \displaystyle\sum _{\substack{1\leq i_{1},\dots,i_{k}\leq n\\1\leq j_{1},\dots,j_{l}\leq n\\ 1\leq r_{1},\dots,r_{k}\leq n\\ 1\leq s_{1},\dots,s_{l}\leq n}}\left(\displaystyle\sum_{a_{1}=1}^{n}R_{i_{1}a_{1}}R_{a_{1}r_{1}}\right)\cdots \left(\displaystyle\sum_{a_{k}=1}^{n}R_{i_{k}a_{k}}R_{a_{k}r_{k}}\right)\left(\displaystyle\sum_{b_{1}=1}^{n}R^{j_{1}b_{1}}R^{b_{1}s_{1}}\right)\cdots\left(\displaystyle\sum_{b_{l}=1}^{n}R^{j_{l}b_{l}}R^{b_{l}s_{l}}\right)F_{j_{1}\dots j_{l}}^{i_{1}\dots i_{k}}F_{s_{1}\dots s_{l}}^{r_{1}\dots r_{k}}\]}

 \[=\displaystyle\sum_{\substack{1\leq a_{1},\dots,a_{k}\leq n \\1\leq b_{1},\dots,b_{l}\leq n}}\displaystyle\sum _{\substack{1\leq i_{1},\dots,i_{k}\leq n\\1\leq j_{1},\dots,j_{l}\leq n\\ 1\leq r_{1},\dots,r_{k}\leq n\\ 1\leq s_{1},\dots,s_{l}\leq n}}R_{i_{1}a_{1}}R_{a_{1}r_{1}}\cdots R_{i_{k}a_{k}}R_{a_{k}r_{k}}R^{j_{1}b_{1}}R^{b_{1}s_{1}}\cdots R^{j_{l}b_{l}}R^{b_{l}s_{l}}F_{j_{1}\dots j_{l}}^{i_{1}\dots i_{k}}F_{s_{1}\dots s_{l}}^{r_{1}\dots r_{k}}\] {\tiny \[=\displaystyle\sum_{\substack{1\leq a_{1},\dots,a_{k}\leq n \\1\leq b_{1},\dots,b_{l}\leq n}}\left(\displaystyle\sum _{\substack{1\leq i_{1},\dots,i_{k}\leq n\\1\leq j_{1},\dots,j_{l}\leq n}}R_{i_{1}a_{1}}\cdots R_{i_{k}a_{k}}R^{j_{1}b_{1}}\cdots R^{j_{l}b_{l}}F_{j_{1}\dots j_{l}}^{i_{1}\dots i_{k}} \right)\left(\displaystyle\sum _{\substack{1\leq r_{1},\dots,r_{k}\leq n\\1\leq s_{1},\dots,s_{l}\leq n}}R_{a_{1}r_{1}}\cdots R_{a_{k}r_{k}}R^{b_{1}s_{1}}\cdots R^{b_{l}s_{l}}F_{s_{1}\dots s_{l}}^{r_{1}\dots r_{k}}\right)\]}
 \[=\displaystyle\sum_{\substack{1\leq a_{1},\dots,a_{k}\leq n \\1\leq b_{1},\dots,b_{l}\leq n}}\left(\displaystyle\sum _{\substack{1\leq i_{1},\dots,i_{k}\leq n\\1\leq j_{1},\dots,j_{l}\leq n}}R_{i_{1}a_{1}}\cdots R_{i_{k}a_{k}}R^{j_{1}b_{1}}\cdots R^{j_{l}b_{l}}F_{j_{1}\dots j_{l}}^{i_{1}\dots i_{k}} \right)^{2}\geq 0,\] which is nonnegative, and where we have used the symmetry of the matrices $R$ and $R^{-1}$. Moreover, if $\langle \mathbf{F},\mathbf{F}\rangle_{\mathbf{g}}=0$, then \[\displaystyle\sum_{\substack{1\leq a_{1},\dots,a_{k}\leq n \\1\leq b_{1},\dots,b_{l}\leq n}}\left(\displaystyle\sum _{\substack{1\leq i_{1},\dots,i_{k}\leq n\\1\leq j_{1},\dots,j_{l}\leq n}}R_{i_{1}a_{1}}\cdots R_{i_{k}a_{k}}R^{j_{1}b_{1}}\cdots R^{j_{l}b_{l}}F_{j_{1}\dots j_{l}}^{i_{1}\dots i_{k}} \right)^{2}=0,\] so for any $a_{1},\dots, a_{k},b_{1},\dots, b_{l}$, \[\displaystyle\sum _{\substack{1\leq i_{1},\dots,i_{k}\leq n\\1\leq j_{1},\dots,j_{l}\leq n}}R_{i_{1}a_{1}}\cdots R_{i_{k}a_{k}}R^{j_{1}b_{1}}\cdots R^{j_{l}b_{l}}F_{j_{1}\dots j_{l}}^{i_{1}\dots i_{k}}=0,\] This can be rewritten as \[\displaystyle\sum_{i_{1}=1}^{n}R_{i_{1}a_{1}}\left(\displaystyle\sum _{\substack{1\leq i_{2},\dots,i_{k}\leq n\\1\leq j_{1},\dots,j_{l}\leq n}}R_{i_{2}a_{2}}\cdots R_{i_{k}a_{k}}R^{j_{1}b_{1}}\cdots R^{j_{l}b_{l}}F_{j_{1}\dots j_{l}}^{i_{1}\dots i_{k}}\right)=0,\] and it holds for every $a_{1}\in\{1,\dots,n\}$. Varying $a_{1}\in\{1,\dots,n\}$ gives a system of equations in $n\times n$ with a unique solution because $R$ is invertible, namely \[\displaystyle\sum _{\substack{1\leq i_{2},\dots,i_{k}\leq n\\1\leq j_{1},\dots,j_{l}\leq n}}R_{i_{2}a_{2}}\cdots R_{i_{k}a_{k}}R^{j_{1}b_{1}}\cdots R^{j_{l}b_{l}}F_{j_{1}\dots j_{l}}^{i_{1}\dots i_{k}}=0\] for every $i_{1}\in\{1,\dots,n\}$.

 Using the invertibility of $R$ and $R^{-1}$ and proceeding in the same way iteratively and inductively from $i_{2}$ through $j_{l-1}$, we conclude that \[\displaystyle\sum_{j_{l}=1}^{n}R^{j_{l}b_l}F_{j_{1}\dots j_{l}}^{i_{1}\dots i_{k}}=0\] for any $i_{1},\dots,i_{k},j_{1},\dots ,j_{l-1}\in\{1,\dots,n\}$. Repeating the argument and using the invertibility of $R^{-1}$, we conclude that, for any $i_{1},\dots,i_{k},j_{1},\dots ,j_{l}\in\{1,\dots,n\}$, \[F_{j_{1}\dots j_{l}}^{i_{1}\dots i_{k}}=0.\] Therefore $\mathbf{F}=0$ and $\langle\text{ },\text{ }\rangle_{\mathbf{g}}$ is positive definite. Its coefficients in any tensor frame are smooth. In an orthonormal frame, the formula reduces to products of Kronecker deltas and shows that the tensor frame in the statement is orthonormal. Finally, decomposable tensors span each fiber, so the required property determines the inner product uniquely.
\end{proof}
 We prove a property that will be useful later.
 \begin{proposition}\label{norma producto tensorial}
 Let $(M,\mathbf{g})$ be a smooth Riemannian manifold with or without boundary.
 If $\mathbf{F}\in \Gamma(T^{(k_{1},l_{1})}(TM))$ and
 $\mathbf{G}\in \Gamma(T^{(k_{2},l_{2})}(TM))$, then
 \[|\mathbf{F}\otimes \mathbf{G}|_{\mathbf{g}}=|\mathbf{F}|_{\mathbf{g}}|\mathbf{G}|_{\mathbf{g}}.\]
 \end{proposition}
\begin{proof}
 Let $(\mathbf{E}_{1},\dots,\mathbf{E}_{n})$ be a smooth local frame defined on an open set $U\subseteq M$. Suppose that $F_{j_{1}\dots j_{l_{1}}}^{i_{1}\dots i_{k_{1}}}$ and $G_{\widetilde{j}_{1}\dots \widetilde{j}_{l_{2}}}^{\widetilde{i}_{1}\dots \widetilde{i}_{k_{2}}}$ are the coordinates of $\mathbf{F}$ and $\mathbf{G}$ on $U$ with respect to this smooth local frame. We know that the coordinates of $\mathbf{F}\otimes \mathbf{G}$ are given by \[(\mathbf{F}\otimes \mathbf{G})_{j_{1}\dots j_{l_{1}+l_{2}}}^{i_{1}\dots i_{k_{1}+k_{2}}}=F_{j_{1}\dots j_{l_{1}}}^{i_{1}\dots i_{k_{1}}}G_{j_{l_{1}+1}\dots j_{l_{1}+l_{2}}}^{i_{k_{1}+1}\dots i_{k_{1}+k_{2}}}.\] Then:
 \[|\mathbf{F}\otimes \mathbf{G}|_{\mathbf{g}}^{2}=\] \[\displaystyle\sum_{\substack{1\leq i_{1},\dots,i_{k_{1}+k_{2}}\leq n\\1\leq j_{1},\dots,j_{l_{1}+l_{2}}\leq n\\ 1\leq r_{1},\dots,r_{k_{1}+k_{2}}\leq n\\ 1\leq s_{1},\dots s_{l_{1}+l_{2}}\leq n}}g_{i_{1}r_{1}}\cdots g_{i_{k_{1}}r_{k_{1}}}g_{i_{k_{1}+1}r_{k_{1}+1}}\cdots g_{i_{k_{1}+k_{2}}r_{k_{1}+k_{2}}}g^{j_{1}s_{1}}\cdots g^{j_{l_{1}}s_{l_{1}}}g^{j_{l_{1}+1}s_{l_{1}+1}}\cdots g^{j_{l_{1}+l_{2}}s_{l_{1}+l_{2}}}\] \[
 F_{j_{1}\dots j_{l_{1}}}^{i_{1}\dots i_{k_{1}}}G_{j_{l_{1}+1}\dots j_{l_{1}+l_{2}}}^{i_{k_{1}+1}\dots i_{k_{1}+k_{2}}}F_{s_{1}\dots s_{l_{1}}}^{r_{1}\dots r_{k_{1}}}G_{s_{l_{1}+1}\dots s_{l_{1}+l_{2}}}^{r_{k_{1}+1}\dots r_{k_{1}+k_{2}}}= \left(\displaystyle\sum_{\substack{1\leq i_{1},\dots,i_{k_{1}}\leq n\\1\leq j_{1},\dots,j_{l_{1}}\leq n\\ 1\leq r_{1},\dots,r_{k_{1}}\leq n\\ 1\leq s_{1},\dots s_{l_{1}}\leq n}}g_{i_{1}r_{1}}\cdots g_{i_{k_{1}}r_{k_{1}}}g^{j_{1}s_{1}}\cdots g^{j_{l_{1}}s_{l_{1}}} F_{j_{1}\dots j_{l_{1}}}^{i_{1}\dots i_{k_{1}}}F_{s_{1}\dots s_{l_{1}}}^{r_{1}\dots r_{k_{1}}}\right)\] \[\left(\displaystyle\sum_{\substack{1\leq i_{k_{1}+1},\dots,i_{k_{1}+k_{2}}\leq n\\1\leq j_{l_{1}+1},\dots,j_{l_{1}+l_{2}}\leq n\\ 1\leq r_{k_{1}+1},\dots,r_{k_{1}+k_{2}}\leq n\\ 1\leq s_{l_{1}+1},\dots s_{l_{1}+l_{2}}\leq n}}g_{i_{k_{1}+1}r_{k_{1}+1}}\cdots g_{i_{k_{1}+k_{2}}r_{k_{1}+k_{2}}}g^{j_{l_{1}+1}s_{l_{1}+1}}\cdots g^{j_{l_{1}+l_{2}}s_{l_{1}+l_{2}}} G_{j_{l_{1}+1}\dots j_{l_{1}+l_{2}}}^{i_{k_{1}+1}\dots i_{k_{1}+k_{2}}}G_{s_{l_{1}+1}\dots s_{l_{1}+l_{2}}}^{r_{k_{1}+1}\dots r_{k_{1}+k_{2}}}\right)\] \[=|\mathbf{F}|_{\mathbf{g}}^{2}|\mathbf{G}|_{\mathbf{g}}^{2}.\]
\end{proof}
\section{Connections}

A vector field takes values in different tangent spaces, so directly subtracting its values at two points has no intrinsic meaning. A connection resolves this problem by specifying how to differentiate a field in a direction and how to compare vectors along a curve. We begin with the Euclidean case, where the directional derivative provides the model that will later be transferred to the manifold.
\begin{example}\label{ej:variedades-riemannianas-recordemos-campos-vectoriales-pensar-funciones-suaves}
 Consider the smooth manifold without boundary $\mathbb{R}^n$ and fields $\mathbf{X},\mathbf{Y}\in \mathfrak{X}(\mathbb{R}^{n})$. Recall that vector fields $\mathbf{X}$ and $\mathbf{Y}$ on $\mathbb{R}^{n}$ can be regarded as smooth functions $\mathbf{X},\mathbf{Y}\colon \mathbb{R}^{n}\longrightarrow \mathbb{R}^{n}$. The concept of directional derivative tells us that, for each $p\in \mathbb{R}^{n}$, \[(D_{\mathbf{X}_{p}}(\mathbf{Y}))(p)=\displaystyle\sum_{i=1}^{n}X^{i}(p)\displaystyle\frac{\partial \mathbf{Y}}{\partial x^{i}}(p).\] This directional derivative gives the rate of change of $\mathbf{Y}$ at $p$ in the direction $\mathbf{X}(p)$. Varying $p\in \mathbb{R}^{n}$ yields a smooth vector field whose pointwise action is $(\nabla_{\mathbf{X}}\mathbf{Y})_{p}:=(D_{\mathbf{X}_{p}}\mathbf{Y})(p)=\displaystyle\sum_{i=1}^{n}X^{i}(p)\displaystyle\frac{\partial \mathbf{Y}}{\partial x^{i}}(p)$. This notion of the derivative of one vector field with respect to another relies heavily on the vector space structure of $\mathbb{R}^{n}$, since the definition of directional derivative involves \[(D_{\mathbf{X}_{p}}(\mathbf{Y}))(p):=\displaystyle\lim_{t\to 0}\displaystyle\frac{\mathbf{Y}(p+t\mathbf{X}_{p})-\mathbf{Y}(p)}{t}\] and in general $\mathbf{Y}(p)\in T_{p}\mathbb{R}^{n}$ and $\mathbf{Y}(p+t\mathbf{X}_{p})\in T_{p+t\mathbf{X}_{p}}\mathbb{R}^{n}$, with both tangent spaces canonically identified with $\mathbb{R}^{n}$. Although this cannot be interpreted directly on manifolds, we can abstract the properties of $\nabla_{\mathbf{X}}\mathbf{Y}$ inherited directly from directional derivatives:
 \begin{enumerate}
 \item If $f\in C^{\infty}(\mathbb{R}^{n})$, then $\nabla_{f\mathbf{X}}\mathbf{Y}=f\nabla_{\mathbf{X}}\mathbf{Y}$.
 \item $\nabla_{\mathbf{X}}\mathbf{Y}$ is $\mathbb{R}$-linear in $\mathbf{Y}$.
 \item For each $f\in C^{\infty}(\mathbb{R}^{n})$, \[\nabla_{\mathbf{X}}(f\mathbf{Y})=f\nabla_{\mathbf{X}}\mathbf{Y}+(\mathbf{X}f)\mathbf{Y}.\]
 \end{enumerate}
 The function $\nabla\colon \mathfrak{X}(\mathbb{R}^{n})\times \mathfrak{X}(\mathbb{R}^{n})\longrightarrow \mathfrak{X}(\mathbb{R}^{n})$ defined in this way is called the \textit{trivial connection} or \textit{Euclidean connection}.
\end{example}
The preceding example leads to the following definition, which provides an intrinsic notion of directional derivative on manifolds:

\begin{definition}\label{definicion 1 conexion}\index{connection}\index{derivative!covariant}\glsadd{conexion}\glsadd{derivada-covariante}
 Let $\pi\colon \mathbf{E}\longrightarrow M$ be a smooth vector bundle over a smooth manifold $M$ with or without boundary. A \textbf{connection on
 $\mathbf{E}$} is a map
 \[
 \nabla\colon \mathfrak{X}(M)\times \Gamma(\mathbf{E})
 \longrightarrow \Gamma(\mathbf{E}),
 \qquad
 (\mathbf{X},\mathbf{u})\longmapsto\nabla_{\mathbf{X}}\mathbf{u},
 \]
 satisfying:
 \begin{enumerate}[label=(\alph*)]
 \item $\nabla_{\mathbf{X}}\mathbf{u}$ is $C^{\infty}(M)$-linear in
 $\mathbf{X}$.
 \item $\nabla_{\mathbf{X}}\mathbf{u}$ is $\mathbb{R}$-linear in
 $\mathbf{u}$.
 \item For each $f\in C^{\infty}(M)$,
 \[
 \nabla_{\mathbf{X}}(f\mathbf{u})
 =f\nabla_{\mathbf{X}}\mathbf{u}+(\mathbf{X}f)\mathbf{u}.
 \]
 \end{enumerate}
\end{definition}
\begin{exercise}\label{conexion mas bilineal}
Let $(M,\mathbf{g})$ be a Riemannian manifold with or without boundary, and let
$\nabla$ be a connection on $TM$. Let
$\mathbf{A}\colon \mathfrak{X}(M)\times\mathfrak{X}(M)
\longrightarrow \mathfrak{X}(M)$ be a bilinear map over $C^{\infty}(M)$,
that is, one satisfying
\[
\mathbf{A}(f\mathbf{X},\mathbf{Y})=f\mathbf{A}(\mathbf{X},\mathbf{Y})
\quad\text{and}\quad
\mathbf{A}(\mathbf{X},f\mathbf{Y})=f\mathbf{A}(\mathbf{X},\mathbf{Y}),
\qquad
\forall f\in C^{\infty}(M),\ \forall \mathbf{X},\mathbf{Y}\in\mathfrak{X}(M).
\]
Show that the map
\[
\widetilde{\nabla}_{\mathbf{X}}\mathbf{Y}
=\nabla_{\mathbf{X}}\mathbf{Y}+\mathbf{A}(\mathbf{X},\mathbf{Y}),
\qquad \mathbf{X},\mathbf{Y}\in\mathfrak{X}(M),
\]
also defines a connection on $TM$.

\end{exercise}

\begin{remark}\label{definicion simbolos christoffel}
 If $\nabla$ is a connection on $TM$ and $(U,(x^{i}))$ is a coordinate chart on $M$, then $\nabla$ determines $n^{3}$ smooth functions $\Gamma_{ij}^{k}\colon U\longrightarrow\mathbb{R}$ such that
 \[
 \nabla_{\boldsymbol{\partial}_{i}}\boldsymbol{\partial}_{j}=\Gamma_{ij}^{k} \boldsymbol{\partial}_{k}.
 \]
 These functions completely describe the local behavior of the connection. They are therefore particularly important and are called the \textit{\textbf{Christoffel symbols}}.
\end{remark}
\begin{example}\label{ej:variedades-riemannianas-conexion-euclidiana-definida-ejemplo-anterior-recordemos}
Let $M=\mathbb{R}^{n}$ be regarded as a smooth manifold without boundary, and let $\nabla$ be the Euclidean connection defined in the preceding example.
Recall that, for each $\mathbf{X},\mathbf{Y}\in\mathfrak{X}(\mathbb{R}^{n})$,
\[
(\nabla_{\mathbf{X}}\mathbf{Y})_{p}=(D_{\mathbf{X}_{p}}\mathbf{Y})(p)=\displaystyle\sum_{i=1}^{n}X^{i}(p)\displaystyle\frac{\partial \mathbf{Y}}{\partial x^{i}}(p).
\]
Now consider the canonical coordinate frame $\{\boldsymbol{\partial}_{1},\dots,\boldsymbol{\partial}_{n}\}$, where $\boldsymbol{\partial}_{i}$ denotes the constant vector field whose $i$th component is $1$ and whose other components are $0$.
Since the fields $\boldsymbol{\partial}_{i}$ are constant, their partial derivatives vanish:
\[
\displaystyle\frac{\partial (\boldsymbol{\partial}_{j})}{\partial x^{i}}=0\quad \text{for all } i,j\in\{1,\dots,n\}.
\]
The definition of $\nabla$ then gives
\[
\nabla_{\boldsymbol{\partial}_{i}}\boldsymbol{\partial}_{j}=0\quad \text{for all } i,j\in\{1,\dots,n\}.
\]
Thus the coefficients of the Euclidean connection in Cartesian coordinates satisfy
\[
\Gamma^{k}_{ij}=0,\quad \text{for all } i,j,k\in\{1,\dots,n\}.
\]

\end{example}

\begin{remark}\label{obs:variedades-riemannianas-conexion-lineal-asi-define-homomorfismo-haces}
 A connection on $\mathbf{E}$ is $C^{\infty}(M)$-linear in $\mathbf{X}$,
 so each $\boldsymbol{\sigma}\in \Gamma(\mathbf{E})$ defines a bundle homomorphism
 \[\nabla\boldsymbol{\sigma}\in \operatorname{Hom}_{C^{\infty}(M)}(\mathfrak{X}(M),\Gamma(\mathbf{E})).\]
 By Proposition~\ref{iso Hom-Gamma-Hom}, there is a canonical isomorphism \[\operatorname{Hom}_{C^{\infty}(M)}(\mathfrak{X}(M),\Gamma(\mathbf{E}))\cong \Gamma(\operatorname{Hom}(TM,\mathbf{E}))\] and by Proposition~\ref{prop:Hom-tensor}, we know that $\operatorname{Hom}(TM,\mathbf{E})\cong (TM)^{*}\otimes \mathbf{E}=T^{(0,1)}(TM)\otimes \mathbf{E}$. We therefore conclude that

 \[ \nabla \boldsymbol{\sigma}\in \operatorname{Hom}_{C^{\infty}(M)}(\mathfrak{X}(M),\Gamma(\mathbf{E}))\cong \Gamma(T^{(0,1)}(TM)\otimes \mathbf{E}),\quad \forall\boldsymbol{\sigma} \in \Gamma(\mathbf{E}).\]
\end{remark}
The following additional properties uniquely characterize the Euclidean connection:
\begin{proposition}\label{prop:conexion-euclidiana-propiedades-y-unicidad}
Let $M=\mathbb{R}^{n}$, a smooth manifold without boundary, have its usual inner product $\langle\cdot,\cdot\rangle$, and let $\nabla$ be the Euclidean connection defined, for
$\mathbf{X},\mathbf{Y}\in\mathfrak{X}(\mathbb{R}^{n})$, by
\[
(\nabla_{\mathbf{X}}\mathbf{Y})_{p}=(D_{\mathbf{X}_{p}}\mathbf{Y})(p)=\displaystyle\sum_{i=1}^{n}X^{i}(p) \partial_{i}\mathbf{Y}(p).
\]
Then $\nabla$ satisfies the following properties:
\begin{enumerate}
 \item For any $\mathbf{X},\mathbf{Y},\mathbf{Z}\in\mathfrak{X}(\mathbb{R}^{n})$,
 \[
 \mathbf{X}(\langle \mathbf{Y},\mathbf{Z}\rangle)=\langle\nabla_{\mathbf{X}}\mathbf{Y},\mathbf{Z}\rangle+\langle \mathbf{Y},\nabla_{\mathbf{X}}\mathbf{Z}\rangle.
 \]
 \item For any $\mathbf{X},\mathbf{Y}\in\mathfrak{X}(\mathbb{R}^{n})$,
 \[
 [\mathbf{X},\mathbf{Y}]=\nabla_{\mathbf{X}}\mathbf{Y}-\nabla_{\mathbf{Y}}\mathbf{X},
 \]
 where $[\mathbf{X},\mathbf{Y}](f)=\mathbf{X}(\mathbf{Y}(f))-\mathbf{Y}(\mathbf{X}(f))$ for every $f\in C^{\infty}(\mathbb{R}^{n})$.
\end{enumerate}
Moreover, $\nabla$ is the unique connection on $\mathbb{R}^{n}$ satisfying both properties.
\end{proposition}

\begin{proof}
Let $\mathbf{X},\mathbf{Y},\mathbf{Z}\in\mathfrak{X}(\mathbb{R}^{n})$.

To prove the first property, write
$\mathbf{Y}=(Y^{1},\dots,Y^{n})$ and $\mathbf{Z}=(Z^{1},\dots,Z^{n})$.
Since $\langle \mathbf{Y},\mathbf{Z}\rangle=\displaystyle\sum_{k=1}^{n}Y^{k}Z^{k}$, we have
\begin{align*}
\mathbf{X}(\langle \mathbf{Y},\mathbf{Z}\rangle)
&= \displaystyle\sum_{i,k=1}^{n}X^{i}\partial_{i}(Y^{k}Z^{k}) \\
&= \displaystyle\sum_{i,k=1}^{n}X^{i}\big((\partial_{i}Y^{k})Z^{k}
+Y^{k}(\partial_{i}Z^{k})\big).
\end{align*}
By the definition of $\nabla$, we have $(\nabla_{\mathbf{X}}\mathbf{Y})^{k}=\displaystyle\sum_{i=1}^{n}X^{i}\partial_{i}Y^{k}$ and $(\nabla_{\mathbf{X}}\mathbf{Z})^{k}=\displaystyle\sum_{i=1}^{n}X^{i}\partial_{i}Z^{k}$. Substituting,
\[
\mathbf{X}(\langle \mathbf{Y},\mathbf{Z}\rangle)
=\displaystyle\sum_{k=1}^{n}(\nabla_{\mathbf{X}}\mathbf{Y})^{k}Z^{k}
+\displaystyle\sum_{k=1}^{n}Y^{k}(\nabla_{\mathbf{X}}\mathbf{Z})^{k}
=\langle\nabla_{\mathbf{X}}\mathbf{Y},\mathbf{Z}\rangle
+\langle \mathbf{Y},\nabla_{\mathbf{X}}\mathbf{Z}\rangle.
\]

For the second property, recall that $[\mathbf{X},\mathbf{Y}](f)=\mathbf{X}(\mathbf{Y}(f))-\mathbf{Y}(\mathbf{X}(f))$. Applying this to $x^{k}$,
\begin{align*}
[\mathbf{X},\mathbf{Y}](x^{k})
&=\mathbf{X}(Y^{k})-\mathbf{Y}(X^{k}) \\
&=\displaystyle\sum_{i=1}^{n}\big(X^{i}\partial_{i}Y^{k}
-Y^{i}\partial_{i}X^{k}\big),
\end{align*}
so $[\mathbf{X},\mathbf{Y}]^{k}=\displaystyle\sum_{i=1}^{n}\big(X^{i}\partial_{i}Y^{k}-Y^{i}\partial_{i}X^{k}\big)$. But also $(\nabla_{\mathbf{X}}\mathbf{Y}-\nabla_{\mathbf{Y}}\mathbf{X})^{k}=\displaystyle\sum_{i=1}^{n}\big(X^{i}\partial_{i}Y^{k}-Y^{i}\partial_{i}X^{k}\big)$, and therefore $[\mathbf{X},\mathbf{Y}]=\nabla_{\mathbf{X}}\mathbf{Y}-\nabla_{\mathbf{Y}}\mathbf{X}$.

\medskip

For uniqueness, let $\widetilde{\nabla}$ be another connection satisfying (a) and (b). In the coordinate frame $\{\boldsymbol{\partial}_{1},\dots,\boldsymbol{\partial}_{n}\}$, write $\widetilde{\nabla}_{\boldsymbol{\partial}_{i}}\boldsymbol{\partial}_{j}=\Gamma^{k}_{ij}\boldsymbol{\partial}_{k}$. From (b), since $[\boldsymbol{\partial}_{i},\boldsymbol{\partial}_{j}]=0$, we obtain $\widetilde{\nabla}_{\boldsymbol{\partial}_{i}}\boldsymbol{\partial}_{j}=\widetilde{\nabla}_{\boldsymbol{\partial}_{j}}\boldsymbol{\partial}_{i}$, that is,
\begin{equation}\label{eq:symmetry}
\Gamma^{k}_{ij}=\Gamma^{k}_{ji}.
\end{equation}

By (a), using the fact that $\langle\boldsymbol{\partial}_{j},\boldsymbol{\partial}_{k}\rangle=\delta_{jk}$ is constant,
\begin{equation}\label{eq:metric}
0=\boldsymbol{\partial}_{i}\langle\boldsymbol{\partial}_{j},\boldsymbol{\partial}_{k}\rangle
=\langle\widetilde{\nabla}_{\boldsymbol{\partial}_{i}}\boldsymbol{\partial}_{j},\boldsymbol{\partial}_{k}\rangle
+\langle\boldsymbol{\partial}_{j},\widetilde{\nabla}_{\boldsymbol{\partial}_{i}}\boldsymbol{\partial}_{k}\rangle
=\Gamma^{k}_{ij}+\Gamma^{j}_{ik},
\end{equation}
whence
\begin{equation}\label{eq:antisymmetry}
\Gamma^{k}_{ij}=-\Gamma^{j}_{ik}.
\end{equation}

Let us show that \eqref{eq:symmetry} and \eqref{eq:antisymmetry} imply $\Gamma^{k}_{ij}=0$. If $j=k$, then $\Gamma^{j}_{ij}=-\Gamma^{j}_{ij}$, hence $\Gamma^{j}_{ij}=0$. If $j\neq k$, successively applying \eqref{eq:symmetry} and \eqref{eq:antisymmetry} gives
\[
\Gamma^{k}_{ij}
=-\Gamma^{j}_{ik}
=-\Gamma^{j}_{ki}
=\Gamma^{i}_{kj}
=\Gamma^{i}_{jk}
=-\Gamma^{k}_{ij},
\]
whence $\Gamma^{k}_{ij}=0$.

Consequently, $\widetilde{\nabla}_{\boldsymbol{\partial}_{i}}\partial_{j}=0$, and $\widetilde{\nabla}$ agrees with the Euclidean connection.
\end{proof}

As the preceding example shows, these properties completely characterize the Euclidean connection and relate it to the standard Riemannian metric on $\mathbb{R}^{n}$. We formalize these two properties on an arbitrary Riemannian manifold as follows:
\begin{definition}\label{def conexion compatible y libre de torsion}\index{metric-compatible and torsion-free connection} Let $(M,\mathbf{g})$ be a Riemannian manifold with or without boundary, and let $\nabla$ be a connection on $TM$. We say that:
 \begin{enumerate}[label=(\alph*)]
 \item $\nabla$ is \textbf{compatible with $\mathbf{g}$} if, for any $\mathbf{X},\mathbf{Y},\mathbf{Z}\in \mathfrak{X}(M)$, $\mathbf{X}(\langle \mathbf{Y},\mathbf{Z}\rangle)=\langle \nabla_{\mathbf{X}}\mathbf{Y},\mathbf{Z}\rangle + \langle \mathbf{Y}, \nabla_{\mathbf{X}}\mathbf{Z}\rangle$.
 \item $\nabla$ is \textbf{torsion-free} if, for any $\mathbf{X},\mathbf{Y}\in\mathfrak{X}(M)$, $[\mathbf{X},\mathbf{Y}]=\nabla_{\mathbf{X}}\mathbf{Y}-\nabla_{\mathbf{Y}}\mathbf{X}$, where \[[\mathbf{X},\mathbf{Y}](f)=\mathbf{X}(\mathbf{Y}(f))-\mathbf{Y}(\mathbf{X}(f)).\]
 \end{enumerate}
\end{definition}
These two conditions are independent, as the following example shows:

\begin{example}\label{ej:variedades-riemannianas-conexion-euclidiana-definida-derivadas-direccionales-dos}
Let $\nabla$ be the Euclidean connection on the manifold without boundary $(\mathbb{R}^{n},\overline{\mathbf{g}})$, defined by directional derivatives.

Consider two modifications of $\nabla$ that alter torsion and metric compatibility, respectively. These modifications are obtained by adding a $C^{\infty}(M)$-bilinear term to the connection, as justified by Exercise \ref{conexion mas bilineal}.
\begin{enumerate}
 \item \textit{A metric-compatible connection with torsion.}

 On $\mathbb{R}^{3}$, define, for $\mathbf{X},\mathbf{Y}\in\mathfrak{X}(\mathbb{R}^{3})$,
 \[
 \widetilde{\nabla}_{\mathbf{X}}\mathbf{Y}=\nabla_{\mathbf{X}}\mathbf{Y}+\mathbf{X}\times \mathbf{Y},
 \]
 where $\times$ denotes the usual cross product on $\mathbb{R}^{3}$.
 The additional term $(\mathbf{X},\mathbf{Y})\mapsto \mathbf{X}\times \mathbf{Y}$ is bilinear over $C^{\infty}(\mathbb{R}^{3})$, so $\widetilde{\nabla}$ defines a connection.

 To verify that $\widetilde{\nabla}$ preserves the metric, note that
 \[
 \mathbf{X}\big(\langle \mathbf{Y},\mathbf{Z}\rangle\big)=\langle \nabla_{\mathbf{X}}\mathbf{Y},\mathbf{Z}\rangle+\langle \mathbf{Y},\nabla_{\mathbf{X}}\mathbf{Z}\rangle,
 \]
 and that the cross product satisfies
 \[
 \langle \mathbf{X}\times \mathbf{Y},\mathbf{Z}\rangle=-\langle \mathbf{Y},\mathbf{X}\times \mathbf{Z}\rangle.
 \]
 Thus,
 \[
\langle \widetilde{\nabla}_{\mathbf{X}}\mathbf{Y},\mathbf{Z}\rangle+\langle \mathbf{Y},\widetilde{\nabla}_{\mathbf{X}}\mathbf{Z}\rangle
 =\langle \nabla_{\mathbf{X}}\mathbf{Y},\mathbf{Z}\rangle+\langle \mathbf{Y},\nabla_{\mathbf{X}}\mathbf{Z}\rangle
 +\langle \mathbf{X}\times \mathbf{Y},\mathbf{Z}\rangle+\langle \mathbf{Y},\mathbf{X}\times \mathbf{Z}\rangle
\]
\[
=\langle \nabla_{\mathbf{X}}\mathbf{Y},\mathbf{Z}\rangle+\langle \mathbf{Y},\nabla_{\mathbf{X}}\mathbf{Z}\rangle
 =\mathbf{X}\big(\langle \mathbf{Y},\mathbf{Z}\rangle\big),
\]
 so $\widetilde{\nabla}$ preserves the metric.

 The torsion of $\widetilde{\nabla}$ is given by
 \[
 \begin{aligned}
 T(\mathbf{X},\mathbf{Y})&=\widetilde{\nabla}_{\mathbf{X}}\mathbf{Y}-\widetilde{\nabla}_{\mathbf{Y}}\mathbf{X}-[\mathbf{X},\mathbf{Y}]\\
 &=(\nabla_{\mathbf{X}}\mathbf{Y}-\nabla_{\mathbf{Y}}\mathbf{X}-[\mathbf{X},\mathbf{Y}])+(\mathbf{X}\times \mathbf{Y}-\mathbf{Y}\times \mathbf{X})\\
 &=2\mathbf{X}\times \mathbf{Y},
 \end{aligned}
 \]
 which is nonzero in general.
 Consequently, $\widetilde{\nabla}$ preserves the metric but is not torsion-free.

 \item \textit{A torsion-free connection that does not preserve the metric.}

 Let $\mathbf{V}\in\mathfrak{X}(\mathbb{R}^{n})$ be a smooth vector field, and define
 \[
 \widehat{\nabla}_{\mathbf{X}}\mathbf{Y}=\nabla_{\mathbf{X}}\mathbf{Y}+\langle \mathbf{X},\mathbf{Y}\rangle \mathbf{V}.
 \]
 The additional term $(\mathbf{X},\mathbf{Y})\mapsto\langle \mathbf{X},\mathbf{Y}\rangle \mathbf{V}$ is symmetric in $\mathbf{X}$ and $\mathbf{Y}$ and bilinear over $C^{\infty}(\mathbb{R}^{n})$, so $\widehat{\nabla}$ is also a connection.

 The torsion vanishes, since
 \[
 \begin{aligned}
 \widehat{\nabla}_{\mathbf{X}}\mathbf{Y}-\widehat{\nabla}_{\mathbf{Y}}\mathbf{X}-[\mathbf{X},\mathbf{Y}]
 &=\nabla_{\mathbf{X}}\mathbf{Y}-\nabla_{\mathbf{Y}}\mathbf{X}-[\mathbf{X},\mathbf{Y}]
 +\langle \mathbf{X},\mathbf{Y}\rangle \mathbf{V}-\langle \mathbf{Y},\mathbf{X}\rangle \mathbf{V}\\
 &=0.
 \end{aligned}
 \]
 We now check whether $\widehat{\nabla}$ preserves the metric. We have
 \[
 \begin{aligned}
 \mathbf{X}\big(\langle \mathbf{Y},\mathbf{Z}\rangle\big)&=\langle \nabla_{\mathbf{X}}\mathbf{Y},\mathbf{Z}\rangle+\langle \mathbf{Y},\nabla_{\mathbf{X}}\mathbf{Z}\rangle,\\
 \langle \widehat{\nabla}_{\mathbf{X}}\mathbf{Y},\mathbf{Z}\rangle&=\langle \nabla_{\mathbf{X}}\mathbf{Y},\mathbf{Z}\rangle+\langle \mathbf{X},\mathbf{Y}\rangle\langle \mathbf{V},\mathbf{Z}\rangle,\\
 \langle \mathbf{Y},\widehat{\nabla}_{\mathbf{X}}\mathbf{Z}\rangle&=\langle \mathbf{Y},\nabla_{\mathbf{X}}\mathbf{Z}\rangle+\langle \mathbf{X},\mathbf{Z}\rangle\langle \mathbf{Y},\mathbf{V}\rangle.
 \end{aligned}
 \]
 Substituting,
 \[
 \mathbf{X}\big(\langle \mathbf{Y},\mathbf{Z}\rangle\big)
 -\langle \widehat{\nabla}_{\mathbf{X}}\mathbf{Y},\mathbf{Z}\rangle
 -\langle \mathbf{Y},\widehat{\nabla}_{\mathbf{X}}\mathbf{Z}\rangle
 =-\langle \mathbf{X},\mathbf{Y}\rangle\langle \mathbf{V},\mathbf{Z}\rangle-\langle \mathbf{X},\mathbf{Z}\rangle\langle \mathbf{Y},\mathbf{V}\rangle,
 \]
 which does not vanish in general.
 Therefore $\widehat{\nabla}$ is torsion-free but does not preserve the Euclidean metric.
\end{enumerate}
\end{example}

Like $\mathbb{R}^{n}$, every Riemannian manifold admits a unique metric-compatible, torsion-free connection, known as the Levi--Civita connection. This result is called the \textit{\textbf{fundamental theorem of Riemannian geometry}}.
\begin{theorem}[Fundamental theorem of Riemannian geometry]\label{teo:variedades-riemannianas-teorema-fundamental-de-la-geometria-riemanniana}\index{fundamental theorem of Riemannian geometry} Let $(M,\mathbf{g})$ be a Riemannian manifold with or without boundary. Then there exists a unique connection $\nabla$ on $TM$ that is compatible with $\mathbf{g}$ and torsion-free. It is called the \textbf{Levi--Civita connection of $\mathbf{g}$}.\end{theorem}
\begin{proof}
If a connection is metric-compatible and torsion-free, then
\[
\begin{aligned}
\mathbf X\langle\mathbf Y,\mathbf Z\rangle
&=\langle\nabla_{\mathbf X}\mathbf Y,\mathbf Z\rangle
 +\langle\mathbf Y,\nabla_{\mathbf X}\mathbf Z\rangle,\\
\mathbf Y\langle\mathbf Z,\mathbf X\rangle
&=\langle\nabla_{\mathbf Y}\mathbf Z,\mathbf X\rangle
 +\langle\mathbf Z,\nabla_{\mathbf Y}\mathbf X\rangle,\\
\mathbf Z\langle\mathbf X,\mathbf Y\rangle
&=\langle\nabla_{\mathbf Z}\mathbf X,\mathbf Y\rangle
 +\langle\mathbf X,\nabla_{\mathbf Z}\mathbf Y\rangle.
\end{aligned}
\]
Add the first two identities, subtract the third, and use
\(\nabla_{\mathbf X}\mathbf Y-\nabla_{\mathbf Y}\mathbf X
=[\mathbf X,\mathbf Y]\). This necessarily gives the Koszul formula
\[
\begin{aligned}
2\langle\nabla_{\mathbf X}\mathbf Y,\mathbf Z\rangle
={}&\mathbf X\langle\mathbf Y,\mathbf Z\rangle
+\mathbf Y\langle\mathbf Z,\mathbf X\rangle
-\mathbf Z\langle\mathbf X,\mathbf Y\rangle\\
&-\langle\mathbf Y,[\mathbf X,\mathbf Z]\rangle
-\langle\mathbf Z,[\mathbf Y,\mathbf X]\rangle
+\langle\mathbf X,[\mathbf Z,\mathbf Y]\rangle.
\end{aligned}
\]
Since the metric is nondegenerate, this identity determines
\(\nabla_{\mathbf X}\mathbf Y\) uniquely.

For existence, denote the right-hand side of the preceding formula by
\(\mathcal K(\mathbf X,\mathbf Y,\mathbf Z)\). Substituting the identities
\[
[f\mathbf X,\mathbf Z]
=f[\mathbf X,\mathbf Z]-\mathbf Z(f)\mathbf X,
\qquad
[\mathbf Z,f\mathbf Y]
=f[\mathbf Z,\mathbf Y]+\mathbf Z(f)\mathbf Y
\]
and canceling the terms containing derivatives of
\(f\) shows that
\[
\mathcal K(f\mathbf X,\mathbf Y,\mathbf Z)
=f\mathcal K(\mathbf X,\mathbf Y,\mathbf Z),
\qquad
\mathcal K(\mathbf X,\mathbf Y,f\mathbf Z)
=f\mathcal K(\mathbf X,\mathbf Y,\mathbf Z).
\]
Likewise,
\[
\mathcal K(\mathbf X,f\mathbf Y,\mathbf Z)
=f\mathcal K(\mathbf X,\mathbf Y,\mathbf Z)
+2\mathbf X(f)\langle\mathbf Y,\mathbf Z\rangle.
\]
Therefore, for each \(\mathbf X,\mathbf Y\), the function
\(\mathbf Z\mapsto\frac12\mathcal K(\mathbf X,\mathbf Y,\mathbf Z)\)
is \(C^\infty(M)\)-linear. Nondegeneracy of the metric defines a unique vector field \(\nabla_{\mathbf X}\mathbf Y\) by
\[
2\langle\nabla_{\mathbf X}\mathbf Y,\mathbf Z\rangle
=\mathcal K(\mathbf X,\mathbf Y,\mathbf Z)
\quad\text{for every }\mathbf Z.
\]
The preceding three identities prove the connection rules
\[
\nabla_{f\mathbf X}\mathbf Y=f\nabla_{\mathbf X}\mathbf Y,\qquad
\nabla_{\mathbf X}(f\mathbf Y)
=f\nabla_{\mathbf X}\mathbf Y+\mathbf X(f)\mathbf Y.
\]
Smoothness is checked by evaluating on a local frame, since the resulting components are combinations of smooth functions and their derivatives.

Finally, the Koszul formula itself gives
\[
\mathcal K(\mathbf X,\mathbf Y,\mathbf Z)
+\mathcal K(\mathbf X,\mathbf Z,\mathbf Y)
=2\mathbf X\langle\mathbf Y,\mathbf Z\rangle
\]
and
\[
\mathcal K(\mathbf X,\mathbf Y,\mathbf Z)
-\mathcal K(\mathbf Y,\mathbf X,\mathbf Z)
=2\langle[\mathbf X,\mathbf Y],\mathbf Z\rangle.
\]
The first identity proves metric compatibility; the second, together with nondegeneracy, proves that the torsion vanishes.
\end{proof}

The Levi--Civita connection will be fundamental in nonlinear analysis on manifolds, since most spaces we work with require a metric-compatible connection. The following corollary collects its local expressions and basic properties.
\begin{corollary}\label{expresion local christoffel levi civita}
 Let $(M,\mathbf{g})$ be a Riemannian manifold with or without boundary, and let $\nabla$ be its Levi--Civita connection.
 \begin{enumerate}
 \item For any $\mathbf{X},\mathbf{Y},\mathbf{Z}\in\mathfrak X(M)$, the \textit{\textbf{Koszul formula}} holds:
 \[\langle \nabla _{\mathbf{X}}\mathbf{Y},\mathbf{Z}\rangle=\displaystyle\frac{1}{2}\left(\mathbf{X}\langle \mathbf{Y},\mathbf{Z}\rangle+
 \mathbf{Y}\langle \mathbf{Z},\mathbf{X}\rangle- \mathbf{Z}\langle \mathbf{X},\mathbf{Y}\rangle-\langle \mathbf{Y},[\mathbf{X},\mathbf{Z}]\rangle -\langle \mathbf{Z},[\mathbf{Y},\mathbf{X}]\rangle + \langle \mathbf{X},[\mathbf{Z},\mathbf{Y}]\rangle \right).\]
 \item In any coordinate chart $(U,\phi)$, we have \[\Gamma_{ij}^{k}=\displaystyle\frac{1}{2}g^{kl}\left(\partial_{i} g_{jl}+\partial_{j} g_{il}-\partial_{l}g_{ij}\right).\]
 \item In a local frame $(\mathbf{E}_{i})$ defined on $U\subseteq M$, if $c_{ij}^{k}\colon U\longrightarrow \mathbb{R}$ are the $n^{3}$ smooth functions defined by \[[\mathbf{E}_{i},\mathbf{E}_{j}]=c_{ij}^{k}\mathbf{E}_{k},\] then the coefficients of the Levi--Civita connection with respect to this frame are \[\Gamma_{ij}^{k}=\displaystyle\frac{1}{2}g^{kl}\left(\mathbf{E}_{i}g_{jl}+\mathbf{E}_{j}g_{il}-\mathbf{E}_{l}g_{ij}-g_{jm}c_{il}^{m}-g_{lm}c_{ji}^{m}+g_{im}c_{lj}^{m}\right).\]
 \item If $(\mathbf{E}_{i})$ is orthonormal, the expression reduces to $\Gamma_{ij}^{k}=\frac{1}{2}\left(c_{ij}^{k}
-c_{ik}^{j}-c_{jk}^{i}\right)$.
 \end{enumerate}
\end{corollary}
\begin{proof}
Part (1) is the formula obtained in the proof of the preceding theorem. In a coordinate frame, the brackets
\([\boldsymbol\partial_i,\boldsymbol\partial_j]\) vanish. Taking
\(\mathbf X=\boldsymbol\partial_i\),
\(\mathbf Y=\boldsymbol\partial_j\), and
\(\mathbf Z=\boldsymbol\partial_l\) gives
\[
2g_{kl}\Gamma_{ij}^k
=\partial_i g_{jl}+\partial_j g_{il}-\partial_l g_{ij}.
\]
Contracting with \(g^{lr}\) gives (2).

For an arbitrary frame \((\mathbf E_i)\), substitute
\([\mathbf E_i,\mathbf E_j]=c_{ij}^k\mathbf E_k\) into the Koszul formula. We obtain
\[
2g_{kl}\Gamma_{ij}^k
=\mathbf E_i g_{jl}+\mathbf E_j g_{il}-\mathbf E_l g_{ij}
-g_{jm}c_{il}^m-g_{lm}c_{ji}^m+g_{im}c_{lj}^m,
\]
and contraction with the inverse metric proves (3). If the frame is orthonormal, \(g_{ij}=\delta_{ij}\) and its derivatives vanish; the formula reduces to
\[
\Gamma_{ij}^k
=\frac12\bigl(c_{ij}^k-c_{ik}^j-c_{jk}^i\bigr),
\]
which is (4).
\end{proof}

The following lemma collects the transformation rules for the Christoffel symbols, the components of the metric, and the components of its inverse.
\begin{lemma}\label{lema: transformacion metrica Christoffel}
Let $(M,\mathbf{g})$ be a Riemannian manifold with or without boundary. Let $(U,\phi)$ and $(V,\psi)$ be two charts on \(M\) such that \(U\cap V\neq \varnothing\). Denote by \(g_{ij}\), \(g^{ij}\), and \(\Gamma_{ij}^{k}\) the components of the metric, the inverse metric, and the Christoffel symbols of the Levi--Civita connection in the coordinates \(\phi =(x^1,\dots,x^n)\), and by \(\widetilde g_{ab}\), \(\widetilde g^{ab}\), and \(\widetilde\Gamma_{ab}^{c}\) the corresponding components in the coordinates \(\psi =(y^1,\dots,y^n)\). Then, for every \(p\in U\cap V\), the following formulas hold:

\begin{enumerate}
 \item The metric components transform according to
 \[
 g_{ij}(p)
 =
 \displaystyle\frac{\partial y^{a}}{\partial x^{i}}(p)
 \displaystyle\frac{\partial y^{b}}{\partial x^{j}}(p)\,
 \widetilde g_{ab}(p).
 \]

 \item The inverse metric components transform according to
 \[
 g^{ij}(p)
 =
 \displaystyle\frac{\partial x^{i}}{\partial y^{a}}(p)
 \displaystyle\frac{\partial x^{j}}{\partial y^{b}}(p)\,
 \widetilde g^{ab}(p).
 \]

 \item The Christoffel symbols transform according to
 \[
 \Gamma_{ij}^{k}(p)
 =
 \displaystyle\frac{\partial x^{k}}{\partial y^{c}}(p)
 \displaystyle\frac{\partial y^{a}}{\partial x^{i}}(p)
 \displaystyle\frac{\partial y^{b}}{\partial x^{j}}(p)\,
 \widetilde\Gamma_{ab}^{c}(p)
 +
 \displaystyle\frac{\partial x^{k}}{\partial y^{c}}(p)
 \displaystyle\frac{\partial^{2}y^{c}}{\partial x^{i}\partial x^{j}}(p).
 \]
\end{enumerate}
\end{lemma}

\begin{proof}
Let \(p\in U\cap V\). We begin with the metric. By the chain rule,
\[
\displaystyle\boldsymbol{\partial}_{i}\Biggr|_{p}
=
\displaystyle\frac{\partial y^{a}}{\partial x^{i}}(p)\,
\displaystyle\boldsymbol{\partial}_{a}\Biggr|_{p}.
\]
Therefore,
\[
g_{ij}(p)
=
\mathbf{g}_p\left(\displaystyle\boldsymbol{\partial}_{i}\Biggr|_{p},\displaystyle\boldsymbol{\partial}_{j}\Biggr|_{p}\right)
=
\displaystyle\frac{\partial y^{a}}{\partial x^{i}}(p)
\displaystyle\frac{\partial y^{b}}{\partial x^{j}}(p)\,
\mathbf{g}_p\left(\displaystyle\boldsymbol{\partial}_{a}\Biggr|_{p},\displaystyle\boldsymbol{\partial}_{b}\Biggr|_{p}\right).
\]
Since \(\mathbf{g}_p\!\left(\displaystyle\boldsymbol{\partial}_{a}\big|_{p},\displaystyle\boldsymbol{\partial}_{b}\big|_{p}\right)=\widetilde g_{ab}(p)\), we obtain
\[
g_{ij}(p)
=
\displaystyle\frac{\partial y^{a}}{\partial x^{i}}(p)
\displaystyle\frac{\partial y^{b}}{\partial x^{j}}(p)\,
\widetilde g_{ab}(p).
\]

We turn to the inverse metric. The preceding identity can be written in matrix form as
\[
G(p)=J(p)^{T}\widetilde G(p)\,J(p),
\]
where \(G(p)=\bigl(g_{ij}(p)\bigr)\), \(\widetilde G(p)=\bigl(\widetilde g_{ab}(p)\bigr)\), and \(J(p)=\left(\displaystyle\frac{\partial y^{a}}{\partial x^{i}}(p)\right)\). Taking inverses in this equality gives
\[
G(p)^{-1}=J(p)^{-1}\widetilde G(p)^{-1}\bigl(J(p)^{T}\bigr)^{-1}.
\]
Since
\[
J(p)^{-1}
=
\left(\displaystyle\frac{\partial x^{i}}{\partial y^{a}}(p)\right),
\]
it follows that
\[
g^{ij}(p)
=
\displaystyle\frac{\partial x^{i}}{\partial y^{a}}(p)
\displaystyle\frac{\partial x^{j}}{\partial y^{b}}(p)\,
\widetilde g^{ab}(p).
\]

Finally, we prove the transformation formula for the Christoffel symbols. By definition,
\[
\nabla_{\boldsymbol{\partial}_{i}}\displaystyle\boldsymbol{\partial}_{j}
=
\Gamma_{ij}^{k}\displaystyle\boldsymbol{\partial}_{k}.
\]
Using again the relation between the coordinate vector fields,
\[
\displaystyle\boldsymbol{\partial}_{i}
=
\displaystyle\frac{\partial y^{a}}{\partial x^{i}}
\displaystyle\boldsymbol{\partial}_{a},
\qquad
\displaystyle\boldsymbol{\partial}_{j}
=
\displaystyle\frac{\partial y^{b}}{\partial x^{j}}
\displaystyle\boldsymbol{\partial}_{b},
\]
and the Leibniz rule for the connection, we have
\[
\nabla_{\boldsymbol{\partial}_{i}}\displaystyle\boldsymbol{\partial}_{j}
=
\nabla_{\frac{\partial y^{a}}{\partial x^{i}}\boldsymbol{\partial}_{a}}
\left(
\displaystyle\frac{\partial y^{b}}{\partial x^{j}}\displaystyle\boldsymbol{\partial}_{b}
\right)
=
\displaystyle\frac{\partial y^{a}}{\partial x^{i}}
\displaystyle\frac{\partial y^{b}}{\partial x^{j}}
\nabla_{\boldsymbol{\partial}_{a}}\displaystyle\boldsymbol{\partial}_{b}
+
\displaystyle\frac{\partial y^{a}}{\partial x^{i}}
\displaystyle\boldsymbol{\partial}_{a}\!\left(\displaystyle\frac{\partial y^{c}}{\partial x^{j}}\right)
\displaystyle\boldsymbol{\partial}_{c}.
\]
Since \(\nabla_{\boldsymbol{\partial}_{a}}\displaystyle\boldsymbol{\partial}_{b}=\widetilde\Gamma_{ab}^{c}\displaystyle\boldsymbol{\partial}_{c}\), this gives
\[
\nabla_{\boldsymbol{\partial}_{i}}\displaystyle\boldsymbol{\partial}_{j}
=
\left(
\displaystyle\frac{\partial y^{a}}{\partial x^{i}}
\displaystyle\frac{\partial y^{b}}{\partial x^{j}}
\widetilde\Gamma_{ab}^{c}
+
\displaystyle\frac{\partial^{2}y^{c}}{\partial x^{i}\partial x^{j}}
\right)
\displaystyle\boldsymbol{\partial}_{c}.
\]
On the other hand, \(\displaystyle\boldsymbol{\partial}_{c}=\displaystyle\frac{\partial x^{k}}{\partial y^{c}}\displaystyle\boldsymbol{\partial}_{k}\). Substituting this identity into the preceding expression, we obtain
\[
\nabla_{\boldsymbol{\partial}_{i}}\displaystyle\boldsymbol{\partial}_{j}
=
\left(
\displaystyle\frac{\partial x^{k}}{\partial y^{c}}
\displaystyle\frac{\partial y^{a}}{\partial x^{i}}
\displaystyle\frac{\partial y^{b}}{\partial x^{j}}
\widetilde\Gamma_{ab}^{c}
+
\displaystyle\frac{\partial x^{k}}{\partial y^{c}}
\displaystyle\frac{\partial^{2}y^{c}}{\partial x^{i}\partial x^{j}}
\right)
\displaystyle\boldsymbol{\partial}_{k}.
\]
Comparing with \(\nabla_{\boldsymbol{\partial}_{i}}\displaystyle\boldsymbol{\partial}_{j}=\Gamma_{ij}^{k}\displaystyle\boldsymbol{\partial}_{k}\), we conclude that
\[
\Gamma_{ij}^{k}(p)
=
\displaystyle\frac{\partial x^{k}}{\partial y^{c}}(p)
\displaystyle\frac{\partial y^{a}}{\partial x^{i}}(p)
\displaystyle\frac{\partial y^{b}}{\partial x^{j}}(p)\,
\widetilde\Gamma_{ab}^{c}(p)
+
\displaystyle\frac{\partial x^{k}}{\partial y^{c}}(p)
\displaystyle\frac{\partial^{2}y^{c}}{\partial x^{i}\partial x^{j}}(p).
\]
This completes the proof.
\end{proof}

The connection on $TM$ induces connections on all tensor bundles. The construction begins with the dual and the tensor product.
\begin{proposition}[Induced connections on tensor bundles]
\label{propiedades de la derivada covariante tensorial}
\index{induced connection!on tensor bundles}
 Let $M$ be a smooth manifold with or without boundary, and let $\nabla$ be a connection on $TM$. Then there exists a unique family of connections, also denoted by $\nabla$, on
 \[
 T^{(k,l)}(TM)=TM^{\otimes k}\otimes(T^*M)^{\otimes l},
 \qquad k,l\in\mathbb N_0,
 \]
 with the following properties.
 \begin{enumerate}[label=(\alph*)]
 \item On $TM$ it agrees with the original connection, and on
 $T^{(0,0)}(TM)=M\times\mathbb R$ it is given by
$\nabla_{\mathbf{X}}f=\mathbf{X}(f)$.
 \item The dual connection satisfies
 \begin{equation}
 \label{eq:conexion-dual-tm}
 (\nabla_{\mathbf{X}}\boldsymbol{\omega})(\mathbf{Y})
 =\mathbf{X}\bigl(\boldsymbol{\omega}(\mathbf{Y})\bigr)-\boldsymbol{\omega}(\nabla_{\mathbf{X}}\mathbf{Y}),
 \end{equation}
 and hence the natural pairing is parallel:
 \[
 \mathbf{X}\bigl(\boldsymbol{\omega}(\mathbf{Y})\bigr)
 =(\nabla_{\mathbf{X}}\boldsymbol{\omega})(\mathbf{Y})+\boldsymbol{\omega}(\nabla_{\mathbf{X}}\mathbf{Y}).
 \]
 \item If $\mathbf{S}\in\Gamma(T^{(r,s)}(TM))$ and
 $\mathbf{T}\in\Gamma(T^{(p,q)}(TM))$, then
 \begin{equation}
 \label{eq:leibniz-conexion-tensorial-temprana}
 \nabla_{\mathbf{X}}(\mathbf{S}\otimes \mathbf{T})
 =(\nabla_{\mathbf{X}}\mathbf{S})\otimes \mathbf{T}+\mathbf{S}\otimes(\nabla_{\mathbf{X}}\mathbf{T}).
 \end{equation}
 \item If $\mathbf{T}\in\Gamma(T^{(k,l)}(TM))$, then
 \begin{align}
 &(\nabla_{\mathbf{X}}\mathbf{T})(\omega^1,\dots,\omega^k,\mathbf{Y}_1,\dots,\mathbf{Y}_l)
 \label{eq:formula-intrinseca-conexion-tensorial}\\
 &\quad=\mathbf{X}\bigl(\mathbf{T}(\omega^1,\dots,\omega^k,\mathbf{Y}_1,\dots,\mathbf{Y}_l)\bigr)\notag\\
 &\qquad-\displaystyle\sum_{a=1}^{k}
 \mathbf{T}(\omega^1,\dots,\nabla_{\mathbf{X}}\omega^a,\dots,\omega^k,
 \mathbf{Y}_1,\dots,\mathbf{Y}_l)\notag\\
 &\qquad-\displaystyle\sum_{b=1}^{l}
 \mathbf{T}(\omega^1,\dots,\omega^k,
 \mathbf{Y}_1,\dots,\nabla_{\mathbf{X}}\mathbf{Y}_b,\dots,\mathbf{Y}_l).\notag
 \end{align}
 \item The connection commutes with every contraction between a contravariant and a covariant index. If $C$ denotes one of these contractions,
 \[
 \nabla_{\mathbf{X}}(CT)=C(\nabla_{\mathbf{X}}\mathbf{T}).
 \]
 It also commutes with permutations of factors, symmetrization, and alternation. In particular, for differential forms,
 \[
 \nabla_{\mathbf{X}}(\boldsymbol{\omega}\wedge\eta)
 =(\nabla_{\mathbf{X}}\boldsymbol{\omega})\wedge\eta+
 \boldsymbol{\omega}\wedge(\nabla_{\mathbf{X}}\eta).
 \]
 \end{enumerate}

 If $(\mathbf{E}_1,\dots,\mathbf{E}_n)$ is a local frame with dual coframe
 $(\boldsymbol{\varepsilon}^1,\dots,\boldsymbol{\varepsilon}^n)$ and
 $\nabla_{\mathbf{E}_a}\mathbf{E}_i=\Gamma^s_{ai}\mathbf{E}_s$, then
 \begin{equation}
 \label{eq:componentes-conexion-tensorial-temprana}
 \begin{aligned}
 (\nabla_{\mathbf{E}_a}\mathbf{T})^{i_1\dots i_k}_{j_1\dots j_l}
 ={}&\mathbf{E}_a\bigl(T^{i_1\dots i_k}_{j_1\dots j_l}\bigr)
 +\displaystyle\sum_{r=1}^{k}\Gamma^{i_r}_{as}
 T^{i_1\dots s\dots i_k}_{j_1\dots j_l}\\
 &-\displaystyle\sum_{t=1}^{l}\Gamma^{s}_{a j_t}
 T^{i_1\dots i_k}_{j_1\dots s\dots j_l}.
 \end{aligned}
 \end{equation}
\end{proposition}

\begin{proof}
We first construct the dual connection. For
$\mathbf{X},\mathbf{Y}\in\mathfrak X(M)$ and $\boldsymbol{\omega}\in\Omega^1(M)$, define
$\nabla_{\mathbf{X}}\boldsymbol{\omega}$ by \eqref{eq:conexion-dual-tm}. The right-hand side is linear over $C^\infty(M)$ in $\mathbf{Y}$. If
$f\in C^\infty(M)$, then
\[
\begin{aligned}
\mathbf{X}\bigl(\boldsymbol{\omega}(f\mathbf{Y})\bigr)-\boldsymbol{\omega}(\nabla_{\mathbf{X}}(f\mathbf{Y}))
&=\mathbf{X}(f)\boldsymbol{\omega}(\mathbf{Y})+f\mathbf{X}\bigl(\boldsymbol{\omega}(\mathbf{Y})\bigr)\\
&\quad-\mathbf{X}(f)\boldsymbol{\omega}(\mathbf{Y})-f\boldsymbol{\omega}(\nabla_{\mathbf{X}}\mathbf{Y})\\
&=f\bigl[\mathbf{X}(\boldsymbol{\omega}(\mathbf{Y}))-\boldsymbol{\omega}(\nabla_{\mathbf{X}}\mathbf{Y})\bigr].
\end{aligned}
\]
By the tensor criterion, this determines a $1$-form. It is smooth because its evaluations on the elements of a local frame are smooth functions.
Moreover,
\[
\nabla_{f\mathbf{X}}\boldsymbol{\omega}=f\nabla_{\mathbf{X}}\boldsymbol{\omega},
\qquad
\nabla_{\mathbf{X}}(f\boldsymbol{\omega})=\mathbf{X}(f)\boldsymbol{\omega}+f\nabla_{\mathbf{X}}\boldsymbol{\omega},
\]
as follows by substitution into \eqref{eq:conexion-dual-tm}. Thus it is a connection on $T^*M$. The pairing identity is the defining formula itself and uniquely determines $\nabla_{\mathbf{X}}\boldsymbol{\omega}$, so the dual connection is unique.

Now consider, more generally, two smooth vector bundles
$\mathbf{E},\mathbf{F}\to M$ equipped with connections $\nabla^{\mathbf{E}}$ and $\nabla^{\mathbf{F}}$. On decomposable local sections, we propose
\begin{equation}
\label{eq:construccion-conexion-producto-prueba-temprana}
\nabla^{\mathbf{E}\otimes \mathbf{F}}_{\mathbf{X}}(\sigma\otimes\tau)
=(\nabla^{\mathbf{E}}_{\mathbf{X}}\sigma)\otimes\tau
+\sigma\otimes(\nabla^{\mathbf{F}}_{\mathbf{X}}\tau).
\end{equation}
This rule respects the balancing relation. Indeed,
\[
\begin{aligned}
\nabla_{\mathbf{X}}((f\sigma)\otimes\tau)
={}&\mathbf{X}(f)\sigma\otimes\tau
+f(\nabla^{\mathbf{E}}_{\mathbf{X}}\sigma)\otimes\tau
+f\sigma\otimes\nabla^{\mathbf{F}}_{\mathbf{X}}\tau,
\end{aligned}
\]
whereas
\[
\begin{aligned}
\nabla_{\mathbf{X}}(\sigma\otimes(f\tau))
={}&f(\nabla^{\mathbf{E}}_{\mathbf{X}}\sigma)\otimes\tau
+\mathbf{X}(f)\sigma\otimes\tau
+f\sigma\otimes\nabla^{\mathbf{F}}_{\mathbf{X}}\tau.
\end{aligned}
\]
Thus both representations of the same tensor yield the same result. On a trivializing open set, the isomorphism
\[
\Gamma(\mathbf{E}|_U)\otimes_{C^\infty(U)}\Gamma(\mathbf{F}|_U)
\longrightarrow\Gamma((\mathbf{E}\otimes \mathbf{F})|_U)
\]
allows us to extend \eqref{eq:construccion-conexion-producto-prueba-temprana}
by additivity to every section. This explicitly proves independence of the representation as a sum of decomposable tensors.
In local frames, the formula contains only smooth coefficients and their derivatives, so it produces a smooth section. It also directly gives
\[
\nabla_{f\mathbf{X}}\mathbf{S}=f\nabla_{\mathbf{X}}\mathbf{S},
\qquad
\nabla_{\mathbf{X}}(f\mathbf{S})=\mathbf{X}(f)\mathbf{S}+f\nabla_{\mathbf{X}}\mathbf{S}.
\]
The local definitions agree on overlaps because they are characterized by the same rule on decomposable sections; they therefore glue to a global connection. Since these sections locally generate
$\Gamma(\mathbf{E}\otimes \mathbf{F})$, the rule determines the connection uniquely.

Apply this construction successively to $TM$ and to the dual connection on
$T^*M$. We obtain connections on
$TM^{\otimes k}\otimes(T^*M)^{\otimes l}$. Uniqueness of the product connection shows that these connections are compatible with the canonical associative identifications, so they do not depend on how the tensor product is parenthesized.

To obtain \eqref{eq:formula-intrinseca-conexion-tensorial}, contract
$\mathbf{T}\otimes\omega^1\otimes\cdots\otimes\omega^k\otimes
\mathbf{Y}_1\otimes\cdots\otimes \mathbf{Y}_l$ to a function. The tensor product rule and the parallel pairing identity give
\[
\begin{aligned}
\mathbf{X}\bigl(\mathbf{T}(\omega^1,\dots,\omega^k,\mathbf{Y}_1,\dots,\mathbf{Y}_l)\bigr)
={}&(\nabla_{\mathbf{X}}\mathbf{T})(\omega^1,\dots,\omega^k,\mathbf{Y}_1,\dots,\mathbf{Y}_l)\\
&+\displaystyle\sum_{a=1}^k
\mathbf{T}(\omega^1,\dots,\nabla_{\mathbf{X}}\omega^a,\dots,\omega^k,
\mathbf{Y}_1,\dots,\mathbf{Y}_l)\\
&+\displaystyle\sum_{b=1}^l
\mathbf{T}(\omega^1,\dots,\omega^k,
\mathbf{Y}_1,\dots,\nabla_{\mathbf{X}}\mathbf{Y}_b,\dots,\mathbf{Y}_l).
\end{aligned}
\]
Solving for the first term yields the intrinsic formula, including the negative signs of all the contravariant and covariant arguments on which $\mathbf{T}$ is evaluated.

An elementary contraction pairs a factor $TM$ with a factor
$T^*M$. When a decomposable tensor is differentiated, the two terms in which the connection acts on these factors combine, by
\eqref{eq:conexion-dual-tm}, into the derivative of their pairing. The other factors remain unchanged. By linearity, it follows that
$\nabla_{\mathbf{X}}(CT)=C(\nabla_{\mathbf{X}}\mathbf{T})$ for every contraction $C$.

If $P_\sigma$ permutes factors of the same type, the product formula applied to a decomposable tensor shows term by term that
$\nabla_{\mathbf{X}}P_\sigma=P_\sigma\nabla_{\mathbf{X}}$. Symmetrization and alternation are linear combinations of these operators and therefore also commute with the connection. Since the wedge product is the normalized alternation of the tensor product, we obtain
\[
\nabla_{\mathbf{X}}(\boldsymbol{\omega}\wedge\eta)
=(\nabla_{\mathbf{X}}\boldsymbol{\omega})\wedge\eta+
\boldsymbol{\omega}\wedge(\nabla_{\mathbf{X}}\eta).
\]

Finally, the identity
$\boldsymbol{\varepsilon}^j(\mathbf{E}_i)=\delta_i^j$ and compatibility with the pairing give
\[
0=\mathbf{E}_a(\delta_i^j)
=(\nabla_{\mathbf{E}_a}\boldsymbol{\varepsilon}^j)(\mathbf{E}_i)
+\boldsymbol{\varepsilon}^j(\nabla_{\mathbf{E}_a}\mathbf{E}_i),
\]
and consequently
\[
\nabla_{\mathbf{E}_a}\boldsymbol{\varepsilon}^j=-\Gamma^j_{as}\boldsymbol{\varepsilon}^s.
\]
Differentiating
\[
\mathbf{T}=T^{i_1\dots i_k}_{j_1\dots j_l}\,
\mathbf{E}_{i_1}\otimes\cdots\otimes \mathbf{E}_{i_k}\otimes
\boldsymbol{\varepsilon}^{j_1}\otimes\cdots\otimes\boldsymbol{\varepsilon}^{j_l}
\]
produces a positive sign for each factor $\mathbf{E}_{i_r}$ and a negative sign for each factor $\boldsymbol{\varepsilon}^{j_t}$. Reindexing these terms gives exactly \eqref{eq:componentes-conexion-tensorial-temprana}.

The preceding rules determine the connection on all elements of a local tensor frame, and the Leibniz rule determines its action on coefficients. This also proves uniqueness of the entire family of induced connections.
\end{proof}
With this proposition, compatibility of the metric with the Levi--Civita connection can be rewritten as follows:
\begin{theorem}\label{nabla X g=0, compatibilidad metrica con conexion levi civita}\index{Levi--Civita connection!metric compatibility}
 Let $(M,\mathbf{g})$ be a Riemannian manifold with or without boundary, and let $\nabla$ be its Levi--Civita connection. The Levi--Civita connection is compatible with $\mathbf{g}$ in the sense that $\nabla_{\mathbf{X}}\mathbf{g}=0$ for every $\mathbf{X}\in \mathfrak{X}(M)$.
\end{theorem}
\begin{proof}
 By Definition~\ref{def conexion compatible y libre de torsion}, since $\nabla$ is a connection compatible with $\mathbf{g}$, we know that, for any $\mathbf{X},\mathbf{Y},\mathbf{Z}\in \mathfrak{X}(M)$, \[\nabla_{\mathbf{X}}\langle \mathbf{Y},\mathbf{Z}\rangle=\langle\nabla_{\mathbf{X}}\mathbf{Y},\mathbf{Z}\rangle+\langle \mathbf{Y},\nabla_{\mathbf{X}}\mathbf{Z}\rangle,\] so Proposition~\ref{propiedades de la derivada covariante tensorial} gives \[\nabla_{\mathbf{X}}\mathbf{g}(\mathbf{Y},\mathbf{Z})=\mathbf{X}(\mathbf{g}(\mathbf{Y},\mathbf{Z}))-\mathbf{g}(\nabla_{\mathbf{X}}\mathbf{Y},\mathbf{Z})-\mathbf{g}(\mathbf{Y},\nabla_{\mathbf{X}}\mathbf{Z})\] \[=\nabla_{\mathbf{X}}\langle \mathbf{Y},\mathbf{Z}\rangle-\langle \nabla_{\mathbf{X}}\mathbf{Y},\mathbf{Z}\rangle-\langle \mathbf{Y},\nabla_{\mathbf{X}}\mathbf{Z}\rangle=0.\]
\end{proof}

In local coordinates, the covariant derivative of a tensor field takes the following form:
\begin{lemma}\label{lema:formula-derivada-covariante-coords}
Let $(M,\mathbf{g})$ be a Riemannian manifold with or without boundary, and let $\nabla$ be a connection on $TM$.
Let $(\mathbf{E}_{1},\dots,\mathbf{E}_{n})$ be a smooth local frame of $TM$ over $U\subseteq M$, and let
$(\boldsymbol{\varepsilon}^{1},\dots,\boldsymbol{\varepsilon}^{n})$ be its dual coframe.
If $\mathbf{F}$ is a tensor field of type $(k,l)$ written in this frame as
\[
\mathbf{F}=\displaystyle\sum_{\substack{1\leq i_{1},\dots,i_{k}\leq n \\[2pt] 1\leq j_{1},\dots,j_{l}\leq n}}
F^{i_{1}\dots i_{k}}_{j_{1}\dots j_{l}}\,
\mathbf{E}_{i_{1}}\otimes\cdots\otimes \mathbf{E}_{i_{k}}\otimes
\boldsymbol{\varepsilon}^{j_{1}}\otimes\cdots\otimes \boldsymbol{\varepsilon}^{j_{l}},
\]
then, for each $a\in\{1,\dots,n\}$,
\[
\bigl(\nabla_{\mathbf{E}_{a}}\mathbf{F}\bigr)^{i_{1}\dots i_{k}}_{j_{1}\dots j_{l}}
=\mathbf{E}_{a}\!\left(F^{i_{1}\dots i_{k}}_{j_{1}\dots j_{l}}\right)
+\displaystyle\sum_{r=1}^{k}\Gamma^{i_{r}}_{a s}\,
F^{i_{1}\dots s \dots i_{k}}_{j_{1}\dots j_{l}}
-\displaystyle\sum_{t=1}^{l}\Gamma^{s}_{a j_{t}}\,
F^{i_{1}\dots i_{k}}_{j_{1}\dots s \dots j_{l}},
\]
where $\Gamma^{c}_{ai}$ are the connection coefficients of $\nabla$ relative to the frame
$(\mathbf{E}_{1},\dots,\mathbf{E}_{n})$, that is,
$\nabla_{\mathbf{E}_{a}}\mathbf{E}_{i}=\Gamma^{c}_{ai}\mathbf{E}_{c}$.
\end{lemma}

\begin{proof}
Let $(\mathbf{E}_{1},\dots,\mathbf{E}_{n})$ be a smooth local frame of $TM$ over $U$, and let
$(\boldsymbol{\varepsilon}^{1},\dots,\boldsymbol{\varepsilon}^{n})$ be the dual coframe, so that
$\langle \boldsymbol{\varepsilon}^{j},\mathbf{E}_{i}\rangle=\delta^{j}_{i}$.
Fix $a\in\{1,\dots,n\}$.

Since $\nabla_{\mathbf{E}_{a}}\mathbf{E}_{i}=\Gamma^{s}_{ai}\mathbf{E}_{s}$, the product rule for the pairing
(Proposition~\ref{propiedades de la derivada covariante tensorial}) applied to
$\langle \boldsymbol{\varepsilon}^{j},\mathbf{E}_{i}\rangle$ gives
\[
0=\mathbf{E}_{a}(\delta^{j}_{i})
=\nabla_{\mathbf{E}_{a}}\langle \boldsymbol{\varepsilon}^{j},\mathbf{E}_{i}\rangle
=\langle \nabla_{\mathbf{E}_{a}}\boldsymbol{\varepsilon}^{j},\mathbf{E}_{i}\rangle
+\langle \boldsymbol{\varepsilon}^{j},\nabla_{\mathbf{E}_{a}}\mathbf{E}_{i}\rangle
=\langle \nabla_{\mathbf{E}_{a}}\boldsymbol{\varepsilon}^{j},\mathbf{E}_{i}\rangle+\Gamma^{j}_{ai},
\]
and therefore $\langle \nabla_{\mathbf{E}_{a}}\boldsymbol{\varepsilon}^{j},\mathbf{E}_{i}\rangle=-\Gamma^{j}_{ai}$ for every $i$, which is equivalent to
\[
\nabla_{\mathbf{E}_{a}}\boldsymbol{\varepsilon}^{j}=-\Gamma^{j}_{as}\boldsymbol{\varepsilon}^{s}.
\]

Writing the tensor field $\mathbf{F}$ in the frame as
\[
\mathbf{F}=\displaystyle\sum_{i_{1},\dots,i_{k}=1}^{n}\displaystyle\sum_{j_{1},\dots,j_{l}=1}^{n}
F^{i_{1}\dots i_{k}}_{j_{1}\dots j_{l}}\,
\mathbf{E}_{i_{1}}\otimes\cdots\otimes \mathbf{E}_{i_{k}}\otimes
\boldsymbol{\varepsilon}^{j_{1}}\otimes\cdots\otimes \boldsymbol{\varepsilon}^{j_{l}},
\]
and using the Leibniz rule for tensor products together with $\nabla_{\mathbf{E}_{a}}f=\mathbf{E}_{a}(f)$ yields a linear combination of tensor fields in the frame, with terms arising from differentiating the coefficients and each factor $\mathbf{E}_{i_{r}}$ and $\boldsymbol{\varepsilon}^{j_{t}}$.
Substituting $\nabla_{\mathbf{E}_{a}}\mathbf{E}_{i_{r}}=\Gamma^{s}_{a i_{r}}\mathbf{E}_{s}$ and
$\nabla_{\mathbf{E}_{a}}\boldsymbol{\varepsilon}^{j_{t}}=-\Gamma^{j_{t}}_{as}\boldsymbol{\varepsilon}^{s}$ and rewriting in the basis of tensor fields associated with the frame, we conclude that
\[
\bigl(\nabla_{\mathbf{E}_{a}}\mathbf{F}\bigr)^{i_{1}\dots i_{k}}_{j_{1}\dots j_{l}}
=\mathbf{E}_{a}\!\left(F^{i_{1}\dots i_{k}}_{j_{1}\dots j_{l}}\right)
+\displaystyle\sum_{r=1}^{k}\Gamma^{i_{r}}_{a s}\,
F^{i_{1}\dots s \dots i_{k}}_{j_{1}\dots j_{l}}
-\displaystyle\sum_{t=1}^{l}\Gamma^{s}_{a j_{t}}\,
F^{i_{1}\dots i_{k}}_{j_{1}\dots s \dots j_{l}}.
\]
\end{proof}

With this proposition, Theorem~\ref{nabla X g=0, compatibilidad metrica con conexion levi civita} takes the following form:
\begin{corollary}\label{compatibilidad levi civita metrica en coordenadas}
Let $(M,\mathbf{g})$ be a Riemannian manifold with or without boundary, with Levi--Civita connection $\nabla$. In a local frame $(\mathbf{E}_{1},\dots,\mathbf{E}_{n})$ with dual coframe $(\boldsymbol{\varepsilon}^{1},\dots,\boldsymbol{\varepsilon}^{n})$, we have
\[
\mathbf{E}_{a}(g_{ij})
=\Gamma^{s}_{ai}g_{sj}
+\Gamma^{s}_{aj}g_{is}.
\]
\end{corollary}

\begin{proof}
Applying Theorem~\ref{nabla X g=0, compatibilidad metrica con conexion levi civita} and Lemma~\ref{lema:formula-derivada-covariante-coords} to the $(0,2)$-tensor field $\mathbf{g}$ gives
\[
0=(\nabla_{a}\mathbf{g})_{ij}
=\mathbf{E}_{a}(g_{ij})
-\Gamma^{s}_{ai}g_{sj}
-\Gamma^{s}_{aj}g_{is}.
\]
It follows directly that
\[
\mathbf{E}_{a}(g_{ij})
=\Gamma^{s}_{ai}g_{sj}
+\Gamma^{s}_{aj}g_{is}.
\]
\end{proof}
\begin{corollary}\label{derivada parcial metrica inversa}
Under the same hypotheses as Corollary~\ref{compatibilidad levi civita metrica en coordenadas}, the inverse metric $g^{ij}$ satisfies the relation
\[
\mathbf{E}_{a}(g^{ij})
= - \left( \Gamma^{i}_{as}g^{sj} + \Gamma^{j}_{as}g^{is} \right).
\]
\end{corollary}

\begin{proof}
Start from the identity $g^{ik}g_{kj} = \delta^{i}_{j}$. Applying $\mathbf{E}_{a}$ and using the Leibniz rule gives
\[
\mathbf{E}_{a}(g^{ik})g_{kj} + g^{ik}E_{a}(g_{kj}) = 0.
\]
Substituting the expression for $\mathbf{E}_{a}(g_{kj})$ given by Corollary~\ref{compatibilidad levi civita metrica en coordenadas}, we have
\[
\mathbf{E}_{a}(g^{ik})g_{kj} + g^{ik}\left( \Gamma^{s}_{ak}g_{sj} + \Gamma^{s}_{aj}g_{ks} \right) = 0.
\]
Expanding and using $g^{ik}g_{ks} = \delta^{i}_{s}$ yields
\[
\mathbf{E}_{a}(g^{ik})g_{kj} + g^{ik}\Gamma^{s}_{ak}g_{sj} + \Gamma^{i}_{aj} = 0.
\]
Multiplying on the right by $g^{jm}$ gives
\[
\mathbf{E}_{a}(g^{ik})\delta^{m}_{k} + g^{ik}\Gamma^{s}_{ak}\delta^{m}_{s} + \Gamma^{i}_{aj}g^{jm} = 0.
\]
That is,
\[
\mathbf{E}_{a}(g^{im}) + g^{ik}\Gamma^{m}_{ak} + \Gamma^{i}_{aj}g^{jm} = 0.
\]
Renaming dummy indices and rearranging, we conclude that
\[
\mathbf{E}_{a}(g^{ij})
= - \left( \Gamma^{i}_{as}g^{sj} + \Gamma^{j}_{as}g^{is} \right).
\]
\end{proof}

 The connection on tensor bundles induced by the Levi--Civita connection is compatible with the induced inner product on tensors:
 \begin{proposition}[Compatibility of the Levi--Civita connection with the tensor inner product]\label{producto interno en tensores compatible con levi civita}\index{Levi--Civita connection!compatibility with the tensor inner product} Suppose that $(M,\mathbf{g})$ is a Riemannian manifold with or without boundary. The connection induced on each tensor bundle by the Levi--Civita connection is compatible with the inner product of tensor fields in the following sense:
 \[\mathbf{X}\langle \mathbf{F},\mathbf{G}\rangle=\langle \nabla_{\mathbf{X}}\mathbf{F},\mathbf{G}\rangle+\langle \mathbf{F},\nabla_{\mathbf{X}}\mathbf{G}\rangle\] for any $\mathbf{X}\in \mathfrak{X}(M)$ and any $\mathbf{F},\mathbf{G}\in \Gamma(T^{(k,l)}(TM))$.

 \end{proposition}
\begin{proof}
The induced connection commutes with tensor products and contractions.
It therefore suffices to verify the identity on simple tensors. Let
\[
\mathbf F=\mathbf Y_1\otimes\cdots\otimes\mathbf Y_k
\otimes\boldsymbol\alpha^1\otimes\cdots\otimes\boldsymbol\alpha^l
\]
and use an analogous expression for \(\mathbf G\), with factors
\(\mathbf Z_i\) and \(\boldsymbol\beta^a\). Their inner product is the product of the pairings of their factors:
\[
\langle\mathbf F,\mathbf G\rangle
=\prod_{i=1}^k\langle\mathbf Y_i,\mathbf Z_i\rangle
\prod_{a=1}^l\langle\boldsymbol\alpha^a,\boldsymbol\beta^a\rangle.
\]
The ordinary Leibniz rule expresses
\(\mathbf X\langle\mathbf F,\mathbf G\rangle\) as a sum of terms in which one of these factors is differentiated. Compatibility of the Levi--Civita connection on \(TM\) gives
\[
\mathbf X\langle\mathbf Y_i,\mathbf Z_i\rangle
=\langle\nabla_{\mathbf X}\mathbf Y_i,\mathbf Z_i\rangle
 +\langle\mathbf Y_i,\nabla_{\mathbf X}\mathbf Z_i\rangle.
\]
On the dual, the definition
\[
(\nabla_{\mathbf X}\boldsymbol\alpha)(\mathbf Y)
=\mathbf X(\boldsymbol\alpha(\mathbf Y))
-\boldsymbol\alpha(\nabla_{\mathbf X}\mathbf Y)
\]
also gives $\nabla_{\mathbf X}(\boldsymbol\beta^\sharp)
=(\nabla_{\mathbf X}\boldsymbol\beta)^\sharp$. To see this, differentiate
$\boldsymbol\beta(\mathbf Y)=\langle\boldsymbol\beta^\sharp,\mathbf Y\rangle$
and use metric compatibility on $TM$; the terms containing
$\nabla_{\mathbf X}\mathbf Y$ cancel. Differentiating
$\langle\boldsymbol\alpha,\boldsymbol\beta\rangle
=\boldsymbol\alpha(\boldsymbol\beta^\sharp)$ now proves the same
compatibility identity for covectors.
Grouping the terms in which a factor of \(\mathbf F\) is differentiated and those in which a factor of \(\mathbf G\) is differentiated yields exactly
\[
\mathbf X\langle\mathbf F,\mathbf G\rangle
=\langle\nabla_{\mathbf X}\mathbf F,\mathbf G\rangle
 +\langle\mathbf F,\nabla_{\mathbf X}\mathbf G\rangle.
\]
Bilinearity extends the equality to arbitrary tensor fields.
\end{proof}

In the total covariant derivative and its iterations, each new derivative index precedes the existing covariant indices.
\begin{proposition}[The total covariant derivative]\label{prop:variedades-riemannianas-la-derivada-covariante-total}\index{covariant derivative!total} Let $M$ be a smooth manifold with or without boundary, and let $\nabla$ be a connection on $TM$. For each $\mathbf{F}\in \Gamma(T^{(k,l)}(TM))$, the map \[\nabla \mathbf{F}\colon \underbrace{\Omega^{1}(M)\times \dots \times \Omega^{1}(M)}_{k\text{ times}}\times \underbrace{\mathfrak{X}(M)\times \dots \times \mathfrak{X}(M)}_{(l+1)\text{ times}}\longrightarrow C^{\infty}(M)\] given by \[(\nabla \mathbf{F})(\boldsymbol{\omega}^{1},\dots,\boldsymbol{\omega}^{k},\mathbf{X},\mathbf{Y}_{1},\dots,\mathbf{Y}_{l})=(\nabla_{\mathbf{X}}\mathbf{F})(\boldsymbol{\omega}^{1},\dots,\boldsymbol{\omega}^{k},\mathbf{Y}_{1},\dots,\mathbf{Y}_{l})\] defines a smooth $(k,l+1)$-tensor field on $M$, called the \textbf{total covariant derivative of $\mathbf{F}$}.

\end{proposition}
\begin{proof}
 By Lemma~\ref{lema de caracterización tensorial}, since $\nabla_{\mathbf{X}}\mathbf{F}$ is a $(k,l)$-tensor field, it is $C^{\infty}(M)$-multilinear in its $k+l$ arguments and is linear over $C^{\infty}(M)$ in $\mathbf{X}$ by the definition of a connection. Thus $\nabla \mathbf{F}$ defines a smooth $(k,l+1)$-tensor field on $M$.
\end{proof}
Lemma~\ref{lema:formula-derivada-covariante-coords} can be expressed in terms of the total covariant derivative as follows:
\begin{corollary}\label{cor:derivada-covariante-total}
Let $(M,\mathbf{g})$ be a Riemannian manifold with or without boundary, with a connection $\nabla$, and let
$(\mathbf{E}_{1},\dots,\mathbf{E}_{n})$ be a smooth local frame of $TM$ over $U\subseteq M$, with dual coframe $(\boldsymbol{\varepsilon}^{1},\dots,\boldsymbol{\varepsilon}^{n})$.
If $\mathbf{F}$ is a tensor field of type $(k,l)$, its total covariant derivative $\nabla \mathbf{F}$ takes the following form in this frame:
\[
\nabla \mathbf{F}
=\displaystyle\sum_{a=1}^{n}\displaystyle\sum_{\substack{1\leq i_{1},\dots,i_{k}\leq n\\ 1\leq j_{1},\dots,j_{l}\leq n}}
(\nabla_{a}\mathbf{F})^{i_{1}\dots i_{k}}_{j_{1}\dots j_{l}}
\mathbf{E}_{i_{1}}\otimes\cdots\otimes \mathbf{E}_{i_{k}}
\otimes \boldsymbol{\varepsilon}^{a}\otimes\boldsymbol{\varepsilon}^{j_{1}}\otimes\cdots\otimes \boldsymbol{\varepsilon}^{j_{l}}.
\]
In particular, its components are
\[
(\nabla \mathbf{F})^{i_{1}\dots i_{k}}_{a j_{1}\dots j_{l}}
=(\nabla_{a}\mathbf{F})^{i_{1}\dots i_{k}}_{j_{1}\dots j_{l}}
=\mathbf{E}_{a}\left(F^{i_{1}\dots i_{k}}_{j_{1}\dots j_{l}}\right)
+\displaystyle\sum_{r=1}^{k}\Gamma^{i_{r}}_{a s}F^{i_{1}\dots s \dots i_{k}}_{j_{1}\dots j_{l}}
-\displaystyle\sum_{t=1}^{l}\Gamma^{s}_{a j_{t}}F^{i_{1}\dots i_{k}}_{j_{1}\dots s \dots j_{l}},
\]
where $\nabla_{a}:=\nabla_{\mathbf{E}_{a}}$ and the connection coefficients $\Gamma^{c}_{ab}$ are defined by $\nabla_{\mathbf{E}_{a}}\mathbf{E}_{b}=\Gamma^{c}_{ab}\mathbf{E}_{c}$.

\end{corollary}
\begin{proof}
The definition of the total covariant derivative implies
\[
(\nabla\mathbf F)
(\boldsymbol\varepsilon^{i_1},\dots,\boldsymbol\varepsilon^{i_k},
\mathbf E_a,\mathbf E_{j_1},\dots,\mathbf E_{j_l})
=(\nabla_{\mathbf E_a}\mathbf F)
(\boldsymbol\varepsilon^{i_1},\dots,\boldsymbol\varepsilon^{i_k},
\mathbf E_{j_1},\dots,\mathbf E_{j_l}).
\]
Lemma~\ref{lema:formula-derivada-covariante-coords} computes the right-hand side as
\[
\mathbf E_a\!\left(F^{i_1\dots i_k}_{j_1\dots j_l}\right)
+\displaystyle\sum_{r=1}^k\Gamma^{i_r}_{as}
F^{i_1\dots s\dots i_k}_{j_1\dots j_l}
-\displaystyle\sum_{t=1}^l\Gamma^s_{a j_t}
F^{i_1\dots i_k}_{j_1\dots s\dots j_l}.
\]
These are precisely the components of \(\nabla\mathbf F\) in the indicated frame and coframe; reconstructing the tensor from its components gives the first formula in the statement.
\end{proof}

\begin{remark}\label{obs:variedades-riemannianas-hemos-definido-campo-tensorial-resulta-ser}
Since $\nabla \mathbf{F}$ has been defined for a $(k,l)$-tensor field and is a $(k,l+1)$-tensor field, we can define $\nabla^{2}\mathbf{F}=\nabla(\nabla \mathbf{F})$, which is a $(k,l+2)$-tensor field. In general, we write $\nabla^{s}\mathbf{F}:=\underbrace{\nabla(\nabla(\dots (\nabla \mathbf{F})))}_{s\text{ times}}\in \Gamma(T^{(k,l+s)}(TM))$. For a function $u\in C^{\infty}(M)$, regarded as a tensor of type $(0,0)$, the same construction defines its iterated covariant derivatives $\nabla^s u$.
\end{remark}
Connections on $TM$ and its tensor bundles are related to the wedge product as follows:

\begin{proposition}[Compatibility of the connection with the wedge product]\label{conexion producto cuña}\index{compatibility of the connection with the wedge product}
Let $M$ be a smooth manifold with or without boundary, and let $\nabla$ be a connection on $TM$.
Then, for any $k,l\geq 0$, if $\boldsymbol{\omega} \in \Gamma(\Lambda^{k}T^{*}M)$ and $\boldsymbol{\eta} \in \Gamma(\Lambda^{l}T^{*}M)$, we have
\[
\nabla_{\mathbf{X}}(\boldsymbol{\omega} \wedge \boldsymbol{\eta}) = (\nabla_{\mathbf{X}}\boldsymbol{\omega})\wedge \boldsymbol{\eta} + \boldsymbol{\omega}\wedge (\nabla_{\mathbf{X}}\boldsymbol{\eta})
\qquad\text{for every } \mathbf{X}\in \mathfrak{X}(M).
\]
\end{proposition}

\begin{proof}
Let $\boldsymbol{\omega}\in \Gamma(\Lambda^{k}T^{*}M)$ and $\boldsymbol{\eta}\in \Gamma(\Lambda^{l}T^{*}M)$. By definition,
\[
\boldsymbol{\omega}\wedge\boldsymbol{\eta}=\displaystyle\frac{1}{k!l!}\displaystyle\sum_{\sigma\in S_{k+l}}(\operatorname{sgn}\sigma){}^{\sigma}(\boldsymbol{\omega}\otimes \boldsymbol{\eta}),
\]
where, for $\mathbf{T}\in \Gamma(T^{(0,k+l)}(TM))$ and $\sigma\in S_{k+l}$, we define
\[
({}^{\sigma}\mathbf{T})(v_{1},\dots,v_{k+l})=\mathbf{T}(v_{\sigma(1)},\dots,v_{\sigma(k+l)}).
\]

For $\mathbf{T}\in \Gamma(T^{(0,k+l)}(TM))$, by
\eqref{eq:formula-intrinseca-conexion-tensorial},
\[
(\nabla_{\mathbf{X}}\mathbf{T})(w_{1},\dots,w_{k+l})=\mathbf{X}(\mathbf{T}(w_{1},\dots,w_{k+l}))-\displaystyle\sum_{j=1}^{k+l}\mathbf{T}(w_{1},\dots,\nabla_{\mathbf{X}}w_{j},\dots,w_{k+l}).
\]
Applying this to ${}^{\sigma}\mathbf{T}$ and evaluating on $(v_{1},\dots,v_{k+l})$ gives
\[
(\nabla_{\mathbf{X}}({}^{\sigma}\mathbf{T}))(v_{1},\dots,v_{k+l})
=\mathbf{X}(\mathbf{T}(v_{\sigma(1)},\dots,v_{\sigma(k+l)}))-\displaystyle\sum_{j=1}^{k+l}\mathbf{T}(v_{\sigma(1)},\dots,\nabla_{\mathbf{X}}v_{\sigma(j)},\dots,v_{\sigma(k+l)})
\]
\[
=(\nabla_{\mathbf{X}}\mathbf{T})(v_{\sigma(1)},\dots,v_{\sigma(k+l)})=({}^{\sigma}(\nabla_{\mathbf{X}}\mathbf{T}))(v_{1},\dots,v_{k+l}),
\]
whence $\nabla_{\mathbf{X}}({}^{\sigma}\mathbf{T})={}^{\sigma}(\nabla_{\mathbf{X}}\mathbf{T})$.

Differentiating the expression for $\boldsymbol{\omega}\wedge\boldsymbol{\eta}$ and using linearity,
\[
\nabla_{\mathbf{X}}(\boldsymbol{\omega}\wedge\boldsymbol{\eta})=\displaystyle\frac{1}{k!l!}\displaystyle\sum_{\sigma\in S_{k+l}}(\operatorname{sgn}\sigma)\nabla_{\mathbf{X}}\big({}^{\sigma}(\boldsymbol{\omega}\otimes \boldsymbol{\eta})\big).
\]
From the preceding calculation,
\[
\nabla_{\mathbf{X}}(\boldsymbol{\omega}\wedge\boldsymbol{\eta})=\displaystyle\frac{1}{k!l!}\displaystyle\sum_{\sigma\in S_{k+l}}(\operatorname{sgn}\sigma){}^{\sigma}\big(\nabla_{\mathbf{X}}(\boldsymbol{\omega}\otimes \boldsymbol{\eta})\big).
\]
By \eqref{eq:leibniz-conexion-tensorial-temprana},
\[
\nabla_{\mathbf{X}}(\boldsymbol{\omega}\otimes \boldsymbol{\eta})=(\nabla_{\mathbf{X}}\boldsymbol{\omega})\otimes \boldsymbol{\eta}+\boldsymbol{\omega}\otimes (\nabla_{\mathbf{X}}\boldsymbol{\eta}).
\]
Substituting,
\[
\nabla_{\mathbf{X}}(\boldsymbol{\omega}\wedge\boldsymbol{\eta})
=\displaystyle\frac{1}{k!l!}\displaystyle\sum_{\sigma\in S_{k+l}}(\operatorname{sgn}\sigma){}^{\sigma}\big((\nabla_{\mathbf{X}}\boldsymbol{\omega})\otimes \boldsymbol{\eta}+\boldsymbol{\omega}\otimes (\nabla_{\mathbf{X}}\boldsymbol{\eta})\big)
\]
\[
=\displaystyle\frac{1}{k!l!}\displaystyle\sum_{\sigma\in S_{k+l}}(\operatorname{sgn}\sigma){}^{\sigma}((\nabla_{\mathbf{X}}\boldsymbol{\omega})\otimes \boldsymbol{\eta})+\displaystyle\frac{1}{k!l!}\displaystyle\sum_{\sigma\in S_{k+l}}(\operatorname{sgn}\sigma){}^{\sigma}(\boldsymbol{\omega}\otimes (\nabla_{\mathbf{X}}\boldsymbol{\eta})).
\]
Recognizing the definition of the wedge product in each term, we conclude that
\[
\nabla_{\mathbf{X}}(\boldsymbol{\omega}\wedge\boldsymbol{\eta})=(\nabla_{\mathbf{X}}\boldsymbol{\omega})\wedge \boldsymbol{\eta}+\boldsymbol{\omega}\wedge (\nabla_{\mathbf{X}}\boldsymbol{\eta}).
\]
\end{proof}
Another operation satisfying the Leibniz rule with respect to connections is interior multiplication:
\begin{lemma}\label{leibniz multiplicacion interior}
Let $(M,\mathbf{g})$ be a Riemannian manifold with or without boundary, $\nabla$ a connection on $TM$, and $(\mathbf{E}_{1},\dots,\mathbf{E}_{n})$ a smooth local frame of $TM$ with dual coframe $(\boldsymbol{\varepsilon}^{1},\dots,\boldsymbol{\varepsilon}^{n})$. For every tensor field $\mathbf{F}\in \Gamma(T^{(k,l+1)}(TM))$, we have
\[
\nabla_{a}\bigl(\iota_{\mathbf{E}_{a}}\mathbf{F}\bigr)
=\iota_{\nabla_{a}\mathbf{E}_{a}}\mathbf{F}+\iota_{\mathbf{E}_{a}}(\nabla_{a}\mathbf{F}).
\]

\end{lemma}

\begin{proof}
By Proposition~\ref{prop:multiplicacion-interior-campo}, if $\mathbf{F}$ has components $F^{i_{1}\dots i_{k}}_{j_{1}\dots j_{l+1}}$ in the frames $(\mathbf{E}_{1},\dots,\mathbf{E}_{n})$ and $(\boldsymbol{\varepsilon}^{1},\dots,\boldsymbol{\varepsilon}^{n})$, then the components of $\iota_{\mathbf{E}_{a}}\mathbf{F}$ are $(\iota_{\mathbf{E}_{a}}\mathbf{F})_{j_{1}\dots j_{l}}^{i_{1}\dots i_{k}}=F_{a j_{1}\dots j_{l}}^{i_{1}\dots i_{k}}$.
Applying Lemma~\ref{lema:formula-derivada-covariante-coords} gives
\begin{equation}
(\nabla_{a}(\iota_{\mathbf{E}_{a}}\mathbf{F}))^{i_{1}\dots i_{k}}_{j_{1}\dots j_{l}}
=\mathbf{E}_{a}\left(F^{i_{1}\dots i_{k}}_{a j_{1}\dots j_{l}}\right)
+\displaystyle\sum_{r=1}^{k}\Gamma^{i_{r}}_{a s}
F^{i_{1}\dots s \dots i_{k}}_{a j_{1}\dots j_{l}}
-\displaystyle\sum_{t=1}^{l}\Gamma^{s}_{a j_{t}}
F^{i_{1}\dots i_{k}}_{a j_{1}\dots s \dots j_{l}}.\end{equation}

On the other hand, Remark~\ref{expresion local multiplicacion interior} gives $\nabla_{a}\mathbf{E}_{a}=\Gamma^{s}_{aa}\mathbf{E}_{s}$, so we have
\[
\bigl(\iota_{\nabla_{a}\mathbf{E}_{a}}\mathbf{F}\bigr)^{i_{1}\dots i_{k}}_{j_{1}\dots j_{l}}
=\Gamma^{s}_{aa}F^{i_{1}\dots i_{k}}_{s j_{1}\dots j_{l}}.
\]

Finally, by Lemma~\ref{lema:formula-derivada-covariante-coords},
\[
\bigl(\iota_{\mathbf{E}_{a}}(\nabla_{a}\mathbf{F})\bigr)^{i_{1}\dots i_{k}}_{j_{1}\dots j_{l}}
=(\nabla_{a}\mathbf{F})^{i_{1}\dots i_{k}}_{a j_{1}\dots j_{l}}
\] \[=\mathbf{E}_{a}\left(F^{i_{1}\dots i_{k}}_{a j_{1}\dots j_{l}}\right)
+\displaystyle\sum_{r=1}^{k}\Gamma^{i_{r}}_{a s}
F^{i_{1}\dots s \dots i_{k}}_{a j_{1}\dots j_{l}}
-\Gamma^{s}_{aa}F^{i_{1}\dots i_{k}}_{s j_{1}\dots j_{l}}
-\displaystyle\sum_{t=1}^{l}\Gamma^{s}_{a j_{t}}
F^{i_{1}\dots i_{k}}_{a j_{1}\dots s \dots j_{l}}.\]

Therefore \[(\nabla_{a}(\iota_{\mathbf{E}_{a}}\mathbf{F}))^{i_{1}\dots i_{k}}_{j_{1}\dots j_{l}}=\bigl(\iota_{\nabla_{a}\mathbf{E}_{a}}\mathbf{F}\bigr)^{i_{1}\dots i_{k}}_{j_{1}\dots j_{l}}+\bigl(\iota_{\mathbf{E}_{a}}(\nabla_{a}\mathbf{F})\bigr)^{i_{1}\dots i_{k}}_{j_{1}\dots j_{l}}.\]
\end{proof}
The notion of a connection on tensor bundles can be generalized to the tensor product of two vector bundles as follows:

\begin{proposition}\label{prop:conexion-producto-abstracta}
Let $\mathbf{E},\mathbf{F}\longrightarrow M$ be smooth vector bundles over a smooth manifold $M$ with or without boundary, equipped with connections $\nabla^{\mathbf{E}}$ and $\nabla^{\mathbf{F}}$ in the sense of Definition~\ref{definicion 1 conexion}. Then:

\begin{enumerate}[label=(\alph*)]

\item If $\mathbf{E}$ and $\mathbf{F}$ have connections $\nabla^{\mathbf{E}}$ and $\nabla^{\mathbf{F}}$, the tensor bundle $\mathbf{E}\otimes \mathbf{F}$ admits a unique connection $\nabla^{\mathbf{E}\otimes \mathbf{F}}\colon \mathfrak{X}(M)\times\Gamma(\mathbf{E}\otimes \mathbf{F})\longrightarrow\Gamma(\mathbf{E}\otimes \mathbf{F})$ determined by the rule
\[
\nabla^{\mathbf{E}\otimes \mathbf{F}}_{\mathbf{X}}(\boldsymbol{\sigma}\otimes\boldsymbol{\tau})
=
(\nabla^{\mathbf{E}}_{\mathbf{X}}\boldsymbol{\sigma})\otimes\boldsymbol{\tau}
+
\boldsymbol{\sigma}\otimes(\nabla^{\mathbf{F}}_{\mathbf{X}}\boldsymbol{\tau}),
\qquad \mathbf{X}\in\mathfrak{X}(M),\ \boldsymbol{\sigma}\in\Gamma(\mathbf{E}),\ \boldsymbol{\tau}\in\Gamma(\mathbf{F}),
\]
and extended by bilinearity and the property
\[
\nabla_{\mathbf{X}}(f\mathbf{T})=f\,\nabla_{\mathbf{X}}\mathbf{T}+(\mathbf{X}f)\,\mathbf{T},\qquad f\in C^{\infty}(M).
\]

\item The connection $\nabla^{\mathbf{E}}$ induces a unique connection $\nabla^{\mathbf{E}^{*}}\colon \mathfrak{X}(M)\times\Gamma(\mathbf{E}^{*})\longrightarrow\Gamma(\mathbf{E}^{*})$ characterized by the identity
\[
\mathbf{X}\langle \mathbf{v},\mathbf{u}\rangle
=
\langle\nabla^{\mathbf{E}^{*}}_{\mathbf{X}}\mathbf{v},\mathbf{u}\rangle
+
\langle \mathbf{v},\nabla^{\mathbf{E}}_{\mathbf{X}}\mathbf{u}\rangle,
\qquad
\mathbf{u}\in\Gamma(\mathbf{E}),\ \mathbf{v}\in\Gamma(\mathbf{E}^{*}),\ \mathbf{X}\in\mathfrak{X}(M).
\]

\item The connection $\nabla^{\mathbf{E}}$ induces a connection $\nabla^{\operatorname{End}(\mathbf{E})}$ on $\operatorname{End}(\mathbf{E})\cong \mathbf{E}^{*}\otimes \mathbf{E}$, defined for every $\mathbf{T}\in\Gamma(\operatorname{End}(\mathbf{E}))$ by
\[
(\nabla^{\operatorname{End}(\mathbf{E})}_{\mathbf{X}}\mathbf{T})(\mathbf{u})
=
\nabla^{\mathbf{E}}_{\mathbf{X}}(\mathbf{T}\mathbf{u})-\mathbf{T}(\nabla^{\mathbf{E}}_{\mathbf{X}}\mathbf{u}),
\qquad
\mathbf{u}\in\Gamma(\mathbf{E}),\ \mathbf{X}\in\mathfrak{X}(M).
\]

\item The connections $\nabla^{\mathbf{E}}$ and $\nabla^{\mathbf{F}}$ induce a connection
\(\nabla^{\operatorname{Hom}(\mathbf{E},\mathbf{F})}\) on the bundle
\(\operatorname{Hom}(\mathbf{E},\mathbf{F})\cong \mathbf{E}^{*}\otimes \mathbf{F}\),
determined by the rule
\[
(\nabla^{\operatorname{Hom}(\mathbf{E},\mathbf{F})}_{\mathbf{X}}\boldsymbol{\Phi})(\mathbf{u})
=
\nabla^{\mathbf{F}}_{\mathbf{X}}(\boldsymbol{\Phi}(\mathbf{u}))-\boldsymbol{\Phi}(\nabla^{\mathbf{E}}_{\mathbf{X}}\mathbf{u}),
\qquad
\boldsymbol{\Phi}\in\Gamma(\operatorname{Hom}(\mathbf{E},\mathbf{F})),\ \mathbf{u}\in\Gamma(\mathbf{E}),\ \mathbf{X}\in\mathfrak{X}(M).
\]

Equivalently, in intrinsic terms, if $\boldsymbol{\Phi}\in \Gamma(\mathbf{E}^{*}\otimes \mathbf{F})$, then $\nabla^{\operatorname{Hom}(\mathbf{E},\mathbf{F})}$ is the connection induced by those on $\mathbf{E}^{*}$ and $\mathbf{F}$ through the tensor product.
\end{enumerate}
\end{proposition}

\begin{proof}
For the tensor product, define the connection on decomposable local sections by the formula in part (a). If
$f\in C^\infty(M)$, the two expansions
\[
\begin{aligned}
\nabla_{\mathbf{X}}((f\mathbf{u})\otimes \mathbf{v})
&=\mathbf{X}(f)\mathbf{u}\otimes \mathbf{v}
+f(\nabla_{\mathbf{X}}^{\mathbf{E}}\mathbf{u})\otimes \mathbf{v}
  +f\mathbf{u}\otimes\nabla_{\mathbf{X}}^{\mathbf{F}}\mathbf{v},\\
\nabla_{\mathbf{X}}(\mathbf{u}\otimes(f\mathbf{v}))
&=f(\nabla_{\mathbf{X}}^{\mathbf{E}}\mathbf{u})\otimes \mathbf{v}
+\mathbf{X}(f)\mathbf{u}\otimes \mathbf{v}
  +f\mathbf{u}\otimes\nabla_{\mathbf{X}}^{\mathbf{F}}\mathbf{v}
\end{aligned}
\]
agree. The rule therefore respects the balancing relations and is independent of the representation of a section as a sum of decomposable tensors. In local frames, it produces smooth coefficients and satisfies the two rules for a connection. Uniqueness follows because decomposable sections locally generate $\mathbf{E}\otimes \mathbf{F}$.

For the dual, set
\[
(\nabla_{\mathbf{X}}^{\mathbf{E}^*}\mathbf{v})(\mathbf{u})
:=\mathbf{X}(\mathbf{v}(\mathbf{u}))
-\mathbf{v}(\nabla_{\mathbf{X}}^{\mathbf{E}}\mathbf{u}).
\]
The same calculation
\[
\mathbf{X}(\mathbf{v}(f\mathbf{u}))
-\mathbf{v}(\nabla_{\mathbf{X}}^{\mathbf{E}}(f\mathbf{u}))
=f\bigl[\mathbf{X}(\mathbf{v}(\mathbf{u}))
-\mathbf{v}(\nabla_{\mathbf{X}}^{\mathbf{E}}\mathbf{u})\bigr]
\]
shows that the right-hand side is linear over $C^\infty(M)$ in $\mathbf{u}$ and defines a smooth section of $\mathbf{E}^*$. Linearity over
$C^\infty(M)$ in $\mathbf{X}$ and the rule
$\nabla_{\mathbf{X}}^{\mathbf{E}^*}(f\mathbf{v})=\mathbf{X}(f)\mathbf{v}+f\nabla_{\mathbf{X}}^{\mathbf{E}^*}\mathbf{v}$ are verified directly.
The pairing identity determines this connection uniquely.

The product connection on $\mathbf{E}^*\otimes \mathbf{E}$ induces the connection on
$\operatorname{End}(\mathbf{E})$. Differentiating the identity
$\mathbf{T}\mathbf{u}$ and using the parallelism of the pairing gives
\[
(\nabla_{\mathbf{X}}^{\operatorname{End}(\mathbf{E})}\mathbf{T})(\mathbf{u})
=\nabla_{\mathbf{X}}^{\mathbf{E}}(\mathbf{T}\mathbf{u})
-\mathbf{T}(\nabla_{\mathbf{X}}^{\mathbf{E}}\mathbf{u}).
\]
The same argument on $\mathbf{E}^*\otimes \mathbf{F}$ gives
\[
(\nabla_{\mathbf{X}}^{\operatorname{Hom}(\mathbf{E},\mathbf{F})}\boldsymbol{\Phi})(\mathbf{u})
=\nabla_{\mathbf{X}}^{\mathbf{F}}(\boldsymbol{\Phi}(\mathbf{u}))
-\boldsymbol{\Phi}(\nabla_{\mathbf{X}}^{\mathbf{E}}\mathbf{u}).
\]
These formulas also prove uniqueness in the last two cases.
\end{proof}

A connection can equivalently be defined as an operator between sections of certain vector bundles, as the following proposition shows. In global analysis, most authors use this equivalence in nearly every setting involving connections.
\begin{proposition}\label{conexion equivalente}
Let $M$ be a smooth manifold with or without boundary, and let
$\mathbf{E}\longrightarrow M$ be a smooth vector bundle. Giving a connection on $\mathbf{E}$ in the sense of Definition~\ref{definicion 1 conexion} is equivalent to giving an operator
\[
 \nabla\colon \Gamma(\mathbf{E})\longrightarrow \Gamma\bigl(T^{(0,1)}(TM)\otimes \mathbf{E}\bigr),
\]
that is linear over $\mathbb{R}$ and satisfies the Leibniz rule
\[
 \nabla(f\boldsymbol{\sigma})=df\otimes\boldsymbol{\sigma}+f\nabla\boldsymbol{\sigma},
 \qquad f\in C^{\infty}(M),\boldsymbol{\sigma}\in\Gamma(\mathbf{E}).
\]

The relationship between the two formulations is described using the canonical isomorphism
\[
\Upsilon:\ \Gamma\bigl(T^{(0,1)}(TM)\otimes \mathbf{E}\bigr)
 \xrightarrow{\ \cong\ }
 \operatorname{Hom}_{C^{\infty}(M)}\bigl(\mathfrak{X}(M),\Gamma(\mathbf{E})\bigr).
\]
Under this identification, for $\mathbf{X}\in\mathfrak{X}(M)$ and $\boldsymbol{\sigma}\in\Gamma(\mathbf{E})$, we have
\[
 \Upsilon(\nabla\boldsymbol{\sigma})(\mathbf{X})=\nabla_{\mathbf{X}}\boldsymbol{\sigma}.
\]
\end{proposition}

\begin{proof}
By Proposition~\ref{prop:Hom-tensor}, applied to $TM$ and $\mathbf{E}$, there is a canonical vector bundle isomorphism
\[
\Phi:\ T^{(0,1)}(TM)\otimes \mathbf{E}\ \xrightarrow{\ \cong\ }\ \operatorname{Hom}(TM,\mathbf{E}),
\qquad
\Phi_p(\omega_p\otimes v)(\mathbf{X}_p)=\omega_p(\mathbf{X}_p)v,
\]
for all $p\in M$, $\omega_p\in T_p^*M$, $v\in \mathbf{E}_p$, and $\mathbf{X}_p\in T_pM$.

Passing to sections gives
\[
\Gamma(\Phi)\colon \ \Gamma\bigl(T^{(0,1)}(TM)\otimes \mathbf{E}\bigr)\ \longrightarrow\ \Gamma\bigl(\operatorname{Hom}(TM,\mathbf{E})\bigr),
\qquad
\bigl(\Gamma(\Phi)(\mathbf{s})\bigr)(p)=\Phi_p\bigl(\mathbf{s}(p)\bigr).
\]
Composing with the identification of Proposition~\ref{iso Hom-Gamma-Hom},
\[
\Theta:\ \Gamma\bigl(\operatorname{Hom}(TM,\mathbf{E})\bigr)\ \xrightarrow{\ \cong\ }\ \operatorname{Hom}_{C^{\infty}(M)}\bigl(\mathfrak{X}(M),\Gamma(\mathbf{E})\bigr),
\qquad
\Theta(\mathbf{A})(\mathbf{X})(p)=\mathbf{A}(p)\bigl(\mathbf{X}_p\bigr),
\]
we obtain the isomorphism of $C^{\infty}(M)$-modules
\[
\Upsilon\ :=\ \Theta\circ\Gamma(\Phi):\ \Gamma\bigl(T^{(0,1)}(TM)\otimes \mathbf{E}\bigr)\ \xrightarrow{\ \cong\ }\ \operatorname{Hom}_{C^{\infty}(M)}\bigl(\mathfrak{X}(M),\Gamma(\mathbf{E})\bigr),
\]
given explicitly by
\[
\bigl(\Upsilon(\mathbf{s})\bigr)(\mathbf{X})(p)=\Theta\bigl(\Gamma(\Phi)(\mathbf{s})\bigr)(\mathbf{X})(p)=\Gamma(\Phi)(\mathbf{s})(p)\bigl(\mathbf{X}_p\bigr)=\Phi_p\bigl(\mathbf{s}(p)\bigr)\bigl(\mathbf{X}_p\bigr).
\]
We denote its inverse by $\Upsilon^{-1}$.

($\Rightarrow$) Suppose we have a connection in the sense of Definition~\ref{definicion 1 conexion}. For fixed $\boldsymbol{\sigma}\in\Gamma(\mathbf{E})$, the map
\[
F_{\boldsymbol{\sigma}}\colon \ \mathfrak{X}(M)\longrightarrow \Gamma(\mathbf{E}),\qquad \mathbf{X}\longmapsto \nabla_{\mathbf{X}}\boldsymbol{\sigma},
\]
is $C^{\infty}(M)$-linear in $\mathbf{X}$. Therefore $F_{\boldsymbol{\sigma}}\in \operatorname{Hom}_{C^{\infty}(M)}\bigl(\mathfrak{X}(M),\Gamma(\mathbf{E})\bigr)$ and, by the isomorphism $\Upsilon^{-1}$, there exists a unique section $\mathbf{s}=\nabla\boldsymbol{\sigma}\in\Gamma\bigl(T^{(0,1)}(TM)\otimes \mathbf{E}\bigr)$ such that
\[
\Upsilon(\nabla\boldsymbol{\sigma})(\mathbf{X})\ =\ \nabla_{\mathbf{X}}\boldsymbol{\sigma},
\qquad \forall \mathbf{X}\in\mathfrak{X}(M).
\]
In particular, for all $p\in M$ and $\mathbf{X}_p\in T_pM$,
\[
\Phi_p\bigl((\nabla\boldsymbol{\sigma})(p)\bigr)(\mathbf{X}_p)\ =\ \nabla_{\mathbf{X}}\boldsymbol{\sigma}(p).
\]
Now let $f\in C^{\infty}(M)$. Using the Leibniz rule for the connection,
\[
\Upsilon\bigl(\nabla(f\boldsymbol{\sigma})\bigr)(\mathbf{X})
=\nabla_{\mathbf{X}}(f\boldsymbol{\sigma})
=(\mathbf{X}f)\boldsymbol{\sigma}+f\nabla_{\mathbf{X}}\boldsymbol{\sigma}.
\]
On the other hand,
\[
\Upsilon(df\otimes\boldsymbol{\sigma})(\mathbf{X})(p)
=\Phi_p\bigl((df\otimes\boldsymbol{\sigma})(p)\bigr)(\mathbf{X}_p)
=df_p(\mathbf{X}_p)\boldsymbol{\sigma}(p)
=(\mathbf{X}f)(p)\boldsymbol{\sigma}(p),
\]
and, since $\Upsilon$ is $C^{\infty}(M)$-linear,
\[
\Upsilon(f\nabla\boldsymbol{\sigma})(\mathbf{X})=f\cdot \Upsilon(\nabla\boldsymbol{\sigma})(\mathbf{X})=f\nabla_{\mathbf{X}}\boldsymbol{\sigma}.
\]
Thus,
\[
\Upsilon\bigl(\nabla(f\boldsymbol{\sigma})\bigr)(\mathbf{X})=\Upsilon\bigl(df\otimes\boldsymbol{\sigma}+f\nabla\boldsymbol{\sigma}\bigr)(\mathbf{X}),
\qquad \forall \mathbf{X}\in\mathfrak{X}(M).
\]
Since $\Upsilon$ is an isomorphism, we conclude that
\[
\nabla(f\boldsymbol{\sigma})=df\otimes\boldsymbol{\sigma}+f\nabla\boldsymbol{\sigma}.
\]

($\Leftarrow$) Conversely, let $\nabla\colon \Gamma(\mathbf{E})\longrightarrow\Gamma\bigl(T^{(0,1)}(TM)\otimes \mathbf{E}\bigr)$ be an operator linear over $\mathbb{R}$ satisfying
\[
\nabla(f\boldsymbol{\sigma})=df\otimes\boldsymbol{\sigma}+f\nabla\boldsymbol{\sigma}.
\]
For $\mathbf{X}\in\mathfrak{X}(M)$ and $\boldsymbol{\sigma}\in\Gamma(\mathbf{E})$, define
\[
\nabla_{\mathbf{X}}\boldsymbol{\sigma}\ :=\ \Upsilon(\nabla\boldsymbol{\sigma})(\mathbf{X}),
\quad\text{that is}\quad
\nabla_{\mathbf{X}}\boldsymbol{\sigma}(p)=\Phi_p\bigl((\nabla\boldsymbol{\sigma})(p)\bigr)(\mathbf{X}_p).
\]
Then $\nabla_{\mathbf{X}}$ is linear in $\boldsymbol{\sigma}$. Moreover, using the Leibniz identity on sections and the $C^{\infty}(M)$-linearity of $\Upsilon$,
\[
\nabla_{\mathbf{X}}(f\boldsymbol{\sigma})=\Upsilon\bigl(\nabla(f\boldsymbol{\sigma})\bigr)(\mathbf{X})
=\Upsilon\bigl(df\otimes\boldsymbol{\sigma}+f\nabla\boldsymbol{\sigma}\bigr)(\mathbf{X})
=\Upsilon(df\otimes\boldsymbol{\sigma})(\mathbf{X})+f\Upsilon(\nabla\boldsymbol{\sigma})(\mathbf{X})
=(\mathbf{X}f)\boldsymbol{\sigma}+f\nabla_{\mathbf{X}}\boldsymbol{\sigma}.
\]
Finally, since $\Upsilon(\nabla\boldsymbol{\sigma})\in \operatorname{Hom}_{C^{\infty}(M)}(\mathfrak{X}(M),\Gamma(\mathbf{E}))$ is $C^{\infty}(M)$-linear in $\mathbf{X}$, we have
\[
\nabla_{f\mathbf{X}}\boldsymbol{\sigma}=\Upsilon(\nabla\boldsymbol{\sigma})(f\mathbf{X})=f\Upsilon(\nabla\boldsymbol{\sigma})(\mathbf{X})=f\nabla_{\mathbf{X}}\boldsymbol{\sigma}.
\]
Therefore $(\mathbf{X},\boldsymbol{\sigma})\mapsto\nabla_{\mathbf{X}}\boldsymbol{\sigma}$ satisfies the properties of Definition~\ref{definicion 1 conexion} and defines a connection.

\end{proof}
Since the two approaches are equivalent, we will sometimes use the two descriptions of a connection interchangeably.
\begin{example}[Trivial connection]\label{ejemplo:conexion-trivial}
Let $M$ be a smooth manifold with or without boundary, and let $\underline{\mathbb{R}}^{k}_{M}=M\times\mathbb{R}^k$ be the trivial bundle of rank $k$ over $M$.
The space of sections $\Gamma(\underline{\mathbb{R}}^{k}_{M})$ agrees with the space $C^{\infty}(M,\mathbb{R}^k)$ of smooth functions taking values in $\mathbb{R}^k$, and we have $(\underline{\mathbb{R}}^{k}_{M})_{p}\cong \mathbb{R}^{k}$ (the fiber spaces of the trivial bundle).

Let $(e_1,\dots,e_k)$ be the canonical basis of $\mathbb{R}^k$, and denote the corresponding constant sections by
$\mathbf{e}_i(p):=(p,e_i)$. If
$\boldsymbol{\sigma}=(\sigma_1,\dots,\sigma_k)\in C^{\infty}(M,\mathbb{R}^k)$, define
\[
\nabla\boldsymbol{\sigma}\ :=\ \displaystyle\sum\limits_{i=1}^k d\sigma_i\otimes \mathbf{e}_i\ \in\ \Gamma\bigl(T^*M\otimes\underline{\mathbb{R}}^k_{M}\bigr).
\]
In other words, for each $p\in M$,
\[
(\nabla\boldsymbol{\sigma})(p)\ =\ \displaystyle\sum\limits_{i=1}^k d\sigma_i(p)\otimes \mathbf{e}_i(p)
\ \in\ T_p^*M\otimes\mathbb{R}^k.
\]

If $\mathbf{X}\in\mathfrak{X}(M)$ is a vector field, then
\[
(\nabla\boldsymbol{\sigma})_{p}(\mathbf{X})\ =\ \displaystyle\sum\limits_{i=1}^k d\sigma_i(p)(\mathbf{X}_p)\mathbf{e}_i(p)
=\ (\mathbf{X}\sigma_1(p),\dots,\mathbf{X}\sigma_k(p)),
\]
that is, the result is obtained by differentiating the section $\boldsymbol{\sigma}$ coordinate by coordinate.

The Leibniz rule is verified directly: for $f\in C^{\infty}(M)$,
\[
\nabla(f\boldsymbol{\sigma})=\displaystyle\sum\limits_{i=1}^k d(f \sigma_i)\otimes \mathbf{e}_i
=\displaystyle\sum\limits_{i=1}^k\bigl((df)\sigma_i+fd\sigma_i\bigr)\otimes \mathbf{e}_i
= df\otimes\boldsymbol{\sigma}+f\nabla\boldsymbol{\sigma}.
\]

Therefore $\nabla$ defines a connection on $\underline{\mathbb{R}}^k_M$, called the \emph{trivial connection}.
\end{example}

Connections and bundle metrics must be compatible with tensor operations. The following constructions allow us to differentiate tensors, duals, and homomorphisms, and will be needed to define Sobolev spaces of sections.

\begin{proposition}[Induced connections on tensor products, duals, and homomorphisms]
\label{conexion en tensores, duales y endomorfismos}
\index{induced connection!on tensor products, duals, and homomorphisms}
Let $M$ be a smooth manifold with or without boundary, and let
$\mathbf{E},\mathbf{F}\to M$ be smooth vector bundles with connections $\nabla^{\mathbf{E}}$ and
$\nabla^{\mathbf{F}}$.
\begin{enumerate}
\item There exists a unique connection $\nabla^{\mathbf{E}\otimes \mathbf{F}}$ on $\mathbf{E}\otimes \mathbf{F}$ such that, for any $\mathbf{X}\in\mathfrak X(M)$,
$\mathbf{u}\in\Gamma(\mathbf{E})$, and $\mathbf{v}\in\Gamma(\mathbf{F})$,
\begin{equation}
\label{eq:conexion-producto-tensorial-haces}
\nabla^{\mathbf{E}\otimes \mathbf{F}}_{\mathbf{X}}(\mathbf{u}\otimes \mathbf{v})
=
(\nabla^{\mathbf{E}}_{\mathbf{X}}\mathbf{u})\otimes \mathbf{v}
+
\mathbf{u}\otimes(\nabla^{\mathbf{F}}_{\mathbf{X}}\mathbf{v}).
\end{equation}

\item There exists a unique connection $\nabla^{\mathbf{E}^*}$ on the dual bundle $\mathbf{E}^*$ characterized by
\begin{equation}
\label{eq:conexion-haz-dual-evaluacion}
\mathbf{X}\bigl(\boldsymbol{\alpha}(\mathbf{u})\bigr)
=
(\nabla^{\mathbf{E}^*}_{\mathbf{X}}\boldsymbol{\alpha})(\mathbf{u})
+
\boldsymbol{\alpha}(\nabla^{\mathbf{E}}_{\mathbf{X}}\mathbf{u})
\end{equation}
for any $\mathbf{X}\in\mathfrak X(M)$, $\boldsymbol{\alpha}\in\Gamma(\mathbf{E}^*)$, and
$\mathbf{u}\in\Gamma(\mathbf{E})$.

\item The bundle $\operatorname{Hom}(\mathbf{E},\mathbf{F})\cong \mathbf{E}^*\otimes \mathbf{F}$ has the induced connection determined by
\begin{equation}
\label{eq:conexion-haz-homomorfismos}
(\nabla^{\operatorname{Hom}(\mathbf{E},\mathbf{F})}_{\mathbf{X}}\mathbf{T})(\mathbf{u})
=
\nabla^{\mathbf{F}}_{\mathbf{X}}(\mathbf{T}\mathbf{u})-\mathbf{T}(\nabla^{\mathbf{E}}_{\mathbf{X}}\mathbf{u}).
\end{equation}
In particular, taking $\mathbf{F}=\mathbf{E}$ gives the connection on
$\operatorname{End}(\mathbf{E})$ and
\[
(\nabla^{\operatorname{End}(\mathbf{E})}_{\mathbf{X}}\mathbf{T})(\mathbf{u})
=
[\nabla^{\mathbf{E}}_{\mathbf{X}},\mathbf{T}]\mathbf{u}.
\]
\end{enumerate}
These constructions are compatible with iterated tensor products and with the canonical identifications between tensor bundles.
\end{proposition}

\begin{proof}
We begin with the tensor product. On a trivializing open set $U$, every section of $\mathbf{E}\otimes \mathbf{F}$ is a finite sum of simple tensor sections.
Define the right-hand side of
\eqref{eq:conexion-producto-tensorial-haces} on these sections and extend it by additivity. We must check that the definition respects the balancing relations of the tensor product. For
$f\in C^\infty(U)$,
\[
\begin{split}
&\nabla^{\mathbf{E}\otimes \mathbf{F}}_{\mathbf{X}}((f\mathbf{u})\otimes \mathbf{v})\\
&\quad=
\mathbf{X}(f)\mathbf{u}\otimes \mathbf{v}
+f(\nabla^{\mathbf{E}}_{\mathbf{X}}\mathbf{u})\otimes \mathbf{v}
+f\mathbf{u}\otimes\nabla^{\mathbf{F}}_{\mathbf{X}}\mathbf{v},
\end{split}
\]
whereas
\[
\begin{split}
&\nabla^{\mathbf{E}\otimes \mathbf{F}}_{\mathbf{X}}(\mathbf{u}\otimes(f\mathbf{v}))\\
&\quad=
(\nabla^{\mathbf{E}}_{\mathbf{X}}\mathbf{u})\otimes f\mathbf{v}
+\mathbf{u}\otimes \mathbf{X}(f)\mathbf{v}
+\mathbf{u}\otimes f\nabla^{\mathbf{F}}_{\mathbf{X}}\mathbf{v}.
\end{split}
\]
The two expressions agree after using
$(f\mathbf{u})\otimes \mathbf{v}=\mathbf{u}\otimes(f\mathbf{v})$. The definition is therefore independent of the representation of a section as a sum of simple tensors. The properties
\[
\nabla^{\mathbf{E}\otimes \mathbf{F}}_{f\mathbf{X}}\mathbf{S}=f\nabla^{\mathbf{E}\otimes \mathbf{F}}_{\mathbf{X}}\mathbf{S}
\]
and
\[
\nabla^{\mathbf{E}\otimes \mathbf{F}}_{\mathbf{X}}(f\mathbf{S})
=\mathbf{X}(f)\mathbf{S}+f\nabla^{\mathbf{E}\otimes \mathbf{F}}_{\mathbf{X}}\mathbf{S}
\]
are checked first on simple tensors and then by additivity. This proves existence of the connection. Uniqueness follows because simple tensor sections locally generate the module of sections of
$\mathbf{E}\otimes \mathbf{F}$.

For the dual bundle, define, for $\boldsymbol{\alpha}\in\Gamma(\mathbf{E}^*)$,
\begin{equation}
\label{eq:definicion-conexion-dual-explicita}
(\nabla^{\mathbf{E}^*}_{\mathbf{X}}\boldsymbol{\alpha})(\mathbf{u})
:=
\mathbf{X}\bigl(\boldsymbol{\alpha}(\mathbf{u})\bigr)-\boldsymbol{\alpha}(\nabla^{\mathbf{E}}_{\mathbf{X}}\mathbf{u}).
\end{equation}
The right-hand side is $C^\infty(M)$-linear in $\mathbf{u}$. For
$f\in C^\infty(M)$,
\[
\begin{split}
&\mathbf{X}\bigl(\boldsymbol{\alpha}(f\mathbf{u})\bigr)-\boldsymbol{\alpha}(\nabla^{\mathbf{E}}_{\mathbf{X}}(f\mathbf{u}))\\
&\quad=
\mathbf{X}(f)\boldsymbol{\alpha}(\mathbf{u})+f\mathbf{X}\bigl(\boldsymbol{\alpha}(\mathbf{u})\bigr)
-\mathbf{X}(f)\boldsymbol{\alpha}(\mathbf{u})-f\boldsymbol{\alpha}(\nabla^{\mathbf{E}}_{\mathbf{X}}\mathbf{u})\\
&\quad=
f\left[\mathbf{X}\bigl(\boldsymbol{\alpha}(\mathbf{u})\bigr)-\boldsymbol{\alpha}(\nabla^{\mathbf{E}}_{\mathbf{X}}\mathbf{u})\right].
\end{split}
\]
Thus \eqref{eq:definicion-conexion-dual-explicita} defines a section of
$\mathbf{E}^*$. The $C^\infty(M)$-linearity in $\mathbf{X}$ and the Leibniz rule in $\boldsymbol{\alpha}$ follow directly from the same formula. By construction,
\eqref{eq:conexion-haz-dual-evaluacion} holds, and this identity determines
$\nabla^{\mathbf{E}^*}$ uniquely.

Now identify $\operatorname{Hom}(\mathbf{E},\mathbf{F})$ with $\mathbf{E}^*\otimes \mathbf{F}$ and use the product connection already constructed. If locally
$\mathbf{T}=\displaystyle\sum_{a=1}^{N}\boldsymbol{\alpha}_a\otimes \mathbf{v}_a$, then
\[
\mathbf{T}\mathbf{u}=\displaystyle\sum_{a=1}^{N}\boldsymbol{\alpha}_a(\mathbf{u})\mathbf{v}_a.
\]
The Leibniz rules and
\eqref{eq:conexion-haz-dual-evaluacion} give
\[
\begin{split}
\nabla^{\mathbf{F}}_{\mathbf{X}}(\mathbf{T}\mathbf{u})
={}&
\displaystyle\sum_{a=1}^{N}
(\nabla^{\mathbf{E}^*}_{\mathbf{X}}\boldsymbol{\alpha}_a)(\mathbf{u})\mathbf{v}_a
+
\displaystyle\sum_{a=1}^{N}
\boldsymbol{\alpha}_a(\nabla^{\mathbf{E}}_{\mathbf{X}}\mathbf{u})\mathbf{v}_a\\
&+
\displaystyle\sum_{a=1}^{N}
\boldsymbol{\alpha}_a(\mathbf{u})\nabla^{\mathbf{F}}_{\mathbf{X}}\mathbf{v}_a.
\end{split}
\]
The first and third terms combine to give
$(\nabla^{\operatorname{Hom}(\mathbf{E},\mathbf{F})}_{\mathbf{X}}\mathbf{T})(\mathbf{u})$, while the second is
$\mathbf{T}(\nabla^{\mathbf{E}}_{\mathbf{X}}\mathbf{u})$. This proves
\eqref{eq:conexion-haz-homomorfismos}. The endomorphism case is the specialization $\mathbf{F}=\mathbf{E}$. Finally, iterations are obtained by successively applying the tensor product construction, and their compatibility with the canonical identifications follows from the uniqueness just proved.
\end{proof}

\begin{lemma}\label{lema:formula-derivada-covariante-coords-E}
Let $M$ be a smooth manifold with or without boundary equipped with a connection $\nabla$ on $TM$, and let $\pi_{\mathbf{E}}\colon \mathbf{E}\longrightarrow M$ be a smooth vector bundle of rank $r$ equipped with a connection $\nabla^{\mathbf{E}}$.
Let $(\mathbf{E}_{1},\dots,\mathbf{E}_{n})$ be a local frame of $TM$ over an open set $U\subseteq M$, $(\boldsymbol{\varepsilon}^{1},\dots,\boldsymbol{\varepsilon}^{n})$ its dual coframe, and $(\mathbf{e}_{1},\dots,\mathbf{e}_{r})$ a local frame of $\mathbf{E}\restriction_{U}$.

Let $\mathbf{F}\in \Gamma(T^{(k,l)}(TM)\otimes \mathbf{E})$ be written in these frames as
\[
\mathbf{F}=
\displaystyle\sum_{\substack{1\le i_{1},\dots,i_{k}\le n\\[2pt]
 1\le j_{1},\dots,j_{l}\le n\\[2pt]
 1\le a\le r}}
F^{i_{1}\dots i_{k}\,a}_{j_{1}\dots j_{l}}\,
\mathbf{E}_{i_{1}}\otimes\cdots\otimes \mathbf{E}_{i_{k}}
\otimes \boldsymbol{\varepsilon}^{j_{1}}\otimes\cdots\otimes \boldsymbol{\varepsilon}^{j_{l}}
\otimes \mathbf{e}_{a}.
\]

Then, for each $c\in\{1,\dots,n\}$, the components of the covariant derivative
\[
\nabla_{\mathbf{E}_{c}}\mathbf{F}:=\big(\nabla^{\mathbf{T}^{(k,l)}\otimes \mathbf{E}}\big)_{\mathbf{E}_{c}}\mathbf{F}
\]
are given by
\[
(\nabla_{\mathbf{E}_{c}}\mathbf{F})^{i_{1}\dots i_{k}\,a}_{j_{1}\dots j_{l}}
=
\mathbf{E}_{c}\big(F^{i_{1}\dots i_{k}\,a}_{j_{1}\dots j_{l}}\big)
+
\displaystyle\sum_{r=1}^{k}\Gamma^{i_{r}}_{c s}\,
F^{i_{1}\dots s\dots i_{k}\,a}_{j_{1}\dots j_{l}}
-
\displaystyle\sum_{t=1}^{l}\Gamma^{s}_{c j_{t}}\,
F^{i_{1}\dots i_{k}\,a}_{j_{1}\dots s\dots j_{l}}
+
\displaystyle\sum_{b=1}^{r}(A_{c})^{a}_{b}\,
F^{i_{1}\dots i_{k}\,b}_{j_{1}\dots j_{l}},
\]
where $\nabla_{\mathbf{E}_{c}}\mathbf{E}_{i}=\Gamma^{s}_{c i}\mathbf{E}_{s}$ and $\nabla^{\mathbf{E}}_{\mathbf{E}_{c}}\mathbf{e}_{b}=(A_{c})^{a}_{b}\,\mathbf{e}_{a}$.
\end{lemma}

\begin{proof}
By Proposition~\ref{prop:conexion-producto-abstracta}, the connection $\nabla$ on $TM$ induces a connection on $T^{(k,l)}(TM)$ characterized by the product rule for the contravariant and covariant factors. In the frame
\[
\Big\{
\mathbf{E}_{i_{1}}\otimes\cdots\otimes \mathbf{E}_{i_{k}}
\otimes \boldsymbol{\varepsilon}^{j_{1}}\otimes\cdots\otimes \boldsymbol{\varepsilon}^{j_{l}}
\Big\}
\]
we have, for every $\mathbf{T}\in\Gamma\bigl(T^{(k,l)}(TM)\bigr)$,
\[
(\nabla_{\mathbf{E}_{c}}\mathbf{T})^{i_{1}\dots i_{k}}_{j_{1}\dots j_{l}}
=
\mathbf{E}_{c}\big(T^{i_{1}\dots i_{k}}_{j_{1}\dots j_{l}}\big)
+
\displaystyle\sum_{r=1}^{k}\Gamma^{i_{r}}_{c s}\,
T^{i_{1}\dots s\dots i_{k}}_{j_{1}\dots j_{l}}
-
\displaystyle\sum_{t=1}^{l}\Gamma^{s}_{c j_{t}}\,
T^{i_{1}\dots i_{k}}_{j_{1}\dots s\dots j_{l}}.
\]

Now write $\mathbf{F}=\displaystyle\sum_{a=1}^{r}\mathbf{F}^{(a)}\otimes \mathbf{e}_{a}$, where, for each $a\in\{1,\dots,r\}$, the tensor field $\mathbf{F}^{(a)}\in\Gamma\bigl(T^{(k,l)}(TM)\bigr)$ has components $\big(\mathbf{F}^{(a)}\big)^{i_{1}\dots i_{k}}_{j_{1}\dots j_{l}}=F^{i_{1}\dots i_{k}\,a}_{j_{1}\dots j_{l}}$.

The induced connection on $T^{(k,l)}(TM)\otimes \mathbf{E}$ is determined by
\[
\nabla_{\mathbf{E}_{c}}(\mathbf{S}\otimes\boldsymbol{\sigma})
=
(\nabla_{\mathbf{E}_{c}}\mathbf{S})\otimes\boldsymbol{\sigma}
+
\mathbf{S}\otimes\nabla^{\mathbf{E}}_{\mathbf{E}_{c}}\boldsymbol{\sigma},
\qquad
\mathbf{S}\in\Gamma\bigl(T^{(k,l)}(TM)\bigr),\ \boldsymbol{\sigma}\in\Gamma(\mathbf{E}),
\]
so, for $\mathbf{F}=\displaystyle\sum_{a=1}^{r}\mathbf{F}^{(a)}\otimes \mathbf{e}_{a}$, we obtain
\[
\nabla_{\mathbf{E}_{c}}\mathbf{F}
=
\displaystyle\sum_{a=1}^{r}
\Big(
(\nabla_{\mathbf{E}_{c}}\mathbf{F}^{(a)})\otimes \mathbf{e}_{a}
+
\mathbf{F}^{(a)}\otimes\nabla^{\mathbf{E}}_{\mathbf{E}_{c}}\mathbf{e}_{a}
\Big).
\]

In the frame $(\mathbf{E}_{1},\dots,\mathbf{E}_{n})$, we have $\nabla_{\mathbf{E}_{c}}\mathbf{E}_{i}=\Gamma^{s}_{c i}\mathbf{E}_{s}$, and, by the definition of the coefficients of the connection $\nabla^{\mathbf{E}}$ in the frame $(\mathbf{e}_{1},\dots,\mathbf{e}_{r})$, $\nabla^{\mathbf{E}}_{\mathbf{E}_{c}}\mathbf{e}_{b}=(A_{c})^{a}_{b}\,\mathbf{e}_{a}$.

Applying Lemma~\ref{lema:formula-derivada-covariante-coords} to each $\mathbf{F}^{(a)}$ gives
\[
(\nabla_{\mathbf{E}_{c}}\mathbf{F}^{(a)})^{i_{1}\dots i_{k}}_{j_{1}\dots j_{l}}
=
\mathbf{E}_{c}\big(F^{i_{1}\dots i_{k}\,a}_{j_{1}\dots j_{l}}\big)
+
\displaystyle\sum_{r=1}^{k}\Gamma^{i_{r}}_{c s}\,
F^{i_{1}\dots s\dots i_{k}\,a}_{j_{1}\dots j_{l}}
-
\displaystyle\sum_{t=1}^{l}\Gamma^{s}_{c j_{t}}\,
F^{i_{1}\dots i_{k}\,a}_{j_{1}\dots s\dots j_{l}}.
\]

On the other hand, using the local expression for the connection on $\mathbf{E}$, $\nabla^{\mathbf{E}}_{\mathbf{E}_{c}}\mathbf{e}_{a}=\displaystyle\sum_{b=1}^{r}(A_{c})^{b}_{a}\,\mathbf{e}_{b}$, we obtain
\[
\mathbf{F}^{(a)}\otimes\nabla^{\mathbf{E}}_{\mathbf{E}_{c}}\mathbf{e}_{a}
=
\displaystyle\sum_{b=1}^{r}(A_{c})^{b}_{a}\,
\big(\mathbf{F}^{(a)}\otimes \mathbf{e}_{b}\big).
\]

Each of these terms is expressed in the basis
\[
\Big\{
\mathbf{E}_{i_{1}}\otimes\cdots\otimes \mathbf{E}_{i_{k}}
\otimes \boldsymbol{\varepsilon}^{j_{1}}\otimes\cdots\otimes \boldsymbol{\varepsilon}^{j_{l}}
\otimes \mathbf{e}_{a}
\Big\},
\]
since $\mathbf{F}^{(a)}\otimes \mathbf{e}_{b}=F^{i_{1}\dots i_{k}\,a}_{j_{1}\dots j_{l}}\,\mathbf{E}_{i_{1}}\otimes\cdots\otimes \mathbf{E}_{i_{k}}\otimes \boldsymbol{\varepsilon}^{j_{1}}\otimes\cdots\otimes \boldsymbol{\varepsilon}^{j_{l}}\otimes \mathbf{e}_{b}$.

Consequently, the term $\mathbf{F}^{(a)}\otimes\nabla^{\mathbf{E}}_{\mathbf{E}_{c}}\mathbf{e}_{a}$ contributes the sum
\[
\displaystyle\sum_{b=1}^{r}(A_{c})^{a}_{b}\,
F^{i_{1}\dots i_{k}\,b}_{j_{1}\dots j_{l}}.
\]
to the components with output index $a$.

To obtain the components of $\nabla_{\mathbf{E}_{c}}\mathbf{F}$, we now need only fix $a\in\{1,\dots,r\}$ and add the contributions of the two preceding blocks to the basis tensor with last factor $\mathbf{e}_{a}$:
\begin{itemize}
\item the block $(\nabla_{\mathbf{E}_{c}}\mathbf{F}^{(a)})\otimes \mathbf{e}_{a}$ contributes
\[
\mathbf{E}_{c}\big(F^{i_{1}\dots i_{k}\,a}_{j_{1}\dots j_{l}}\big)
+
\displaystyle\sum_{r=1}^{k}\Gamma^{i_{r}}_{c s}\,
F^{i_{1}\dots s\dots i_{k}\,a}_{j_{1}\dots j_{l}}
-
\displaystyle\sum_{t=1}^{l}\Gamma^{s}_{c j_{t}}\,
F^{i_{1}\dots i_{k}\,a}_{j_{1}\dots s\dots j_{l}},
\]
\item and the block $\mathbf{F}^{(b)}\otimes\nabla^{\mathbf{E}}_{\mathbf{E}_{c}}\mathbf{e}_{b}$ contributes
\[
\displaystyle\sum_{b=1}^{r}(A_{c})^{a}_{b}\,
F^{i_{1}\dots i_{k}\,b}_{j_{1}\dots j_{l}}.
\]
\end{itemize}

Grouping these contributions gives
\[
(\nabla_{\mathbf{E}_{c}}\mathbf{F})^{i_{1}\dots i_{k}\,a}_{j_{1}\dots j_{l}}
=
\mathbf{E}_{c}\big(F^{i_{1}\dots i_{k}\,a}_{j_{1}\dots j_{l}}\big)
+
\displaystyle\sum_{r=1}^{k}\Gamma^{i_{r}}_{c s}\,
F^{i_{1}\dots s\dots i_{k}\,a}_{j_{1}\dots j_{l}}
-
\displaystyle\sum_{t=1}^{l}\Gamma^{s}_{c j_{t}}\,
F^{i_{1}\dots i_{k}\,a}_{j_{1}\dots s\dots j_{l}}
+
\displaystyle\sum_{b=1}^{r}(A_{c})^{a}_{b}\,
F^{i_{1}\dots i_{k}\,b}_{j_{1}\dots j_{l}},
\]
which is the desired formula.
\end{proof}

\begin{proposition}[Formulas for the second and third covariant derivatives]
\label{prop:formulas-segunda-tercera-derivada-covariante}
\index{covariant derivative!second and third}
Let $M$ be a smooth manifold with or without boundary equipped with a connection $\nabla^M$ on $TM$, and let
$\mathbf{E}\to M$ be a smooth vector bundle with a connection $\nabla^{\mathbf{E}}$. On the bundles $T^{(0,q)}(TM)\otimes \mathbf{E}$, use the product connections induced by
$\nabla^M$ and $\nabla^{\mathbf{E}}$. For $\mathbf{u}\in\Gamma(\mathbf{E})$, define
\[
\nabla^0\mathbf{u}:=\mathbf{u},
\qquad
\nabla^{q+1}\mathbf{u}:=\nabla(\nabla^q\mathbf{u}),
\qquad q\in\mathbb N_0.
\]
Thus,
\[
\nabla^q\mathbf{u}\in\Gamma\bigl(T^{(0,q)}(TM)\otimes \mathbf{E}\bigr),
\]
and the new covariant index introduced at each step is placed in the first argument. Equivalently,
\[
(\nabla^{q+1}\mathbf{u})(\mathbf{X}_0,\mathbf{X}_1,\ldots,\mathbf{X}_q)
=\bigl(\nabla_{\mathbf{X}_0}(\nabla^q\mathbf{u})\bigr)(\mathbf{X}_1,\ldots,\mathbf{X}_q).
\]
In particular, for any
$\mathbf{X},\mathbf{Y}\in\mathfrak X(M)$,
\begin{equation}
\label{eq:formula-segunda-derivada-covariante-seccion}
(\nabla^2\mathbf{u})(\mathbf{X},\mathbf{Y})
=
\nabla^{\mathbf{E}}_{\mathbf{X}}\nabla^{\mathbf{E}}_{\mathbf{Y}}\mathbf{u}
-
\nabla^{\mathbf{E}}_{\nabla^M_{\mathbf{X}}\mathbf{Y}}\mathbf{u}.
\end{equation}
Moreover, for any $\mathbf{X},\mathbf{Y},\mathbf{Z}\in\mathfrak X(M)$,
\begin{equation}
\label{eq:formula-tercera-derivada-covariante-seccion}
\begin{aligned}
(\nabla^3\mathbf{u})(\mathbf{X},\mathbf{Y},\mathbf{Z})={}&
\nabla^{\mathbf{E}}_{\mathbf{X}}\nabla^{\mathbf{E}}_{\mathbf{Y}}\nabla^{\mathbf{E}}_{\mathbf{Z}}\mathbf{u}-\nabla^{\mathbf{E}}_{\mathbf{X}}\nabla^{\mathbf{E}}_{\nabla^M_{\mathbf{Y}}\mathbf{Z}}\mathbf{u}
-\nabla^{\mathbf{E}}_{\nabla^M_{\mathbf{X}}\mathbf{Y}}\nabla^{\mathbf{E}}_{\mathbf{Z}}\mathbf{u}\\
&+\nabla^{\mathbf{E}}_{\nabla^M_{\nabla^M_{\mathbf{X}}\mathbf{Y}}\mathbf{Z}}\mathbf{u}
-\nabla^{\mathbf{E}}_{\mathbf{Y}}\nabla^{\mathbf{E}}_{\nabla^M_{\mathbf{X}}\mathbf{Z}}\mathbf{u}
+\nabla^{\mathbf{E}}_{\nabla^M_{\mathbf{Y}}\nabla^M_{\mathbf{X}}\mathbf{Z}}\mathbf{u}.
\end{aligned}
\end{equation}
\end{proposition}

\begin{proof}
The first covariant derivative is interpreted as the $1$-form with values in $\mathbf{E}$ given by
\[
(\nabla \mathbf{u})(\mathbf{X})=\nabla^{\mathbf{E}}_{\mathbf{X}}\mathbf{u}.
\]
The product connection on $T^*M\otimes \mathbf{E}$ satisfies the Leibniz rule for evaluation. Hence,
\[
\begin{aligned}
(\nabla^2\mathbf{u})(\mathbf{X},\mathbf{Y})
&=\bigl(\nabla_{\mathbf{X}}(\nabla \mathbf{u})\bigr)(\mathbf{Y})
=\nabla^{\mathbf{E}}_{\mathbf{X}}\bigl((\nabla \mathbf{u})(\mathbf{Y})\bigr)-(\nabla \mathbf{u})(\nabla^M_{\mathbf{X}}\mathbf{Y})\\
&=\nabla^{\mathbf{E}}_{\mathbf{X}}\nabla^{\mathbf{E}}_{\mathbf{Y}}\mathbf{u}-\nabla^{\mathbf{E}}_{\nabla^M_{\mathbf{X}}\mathbf{Y}}\mathbf{u}.
\end{aligned}
\]
which proves
\eqref{eq:formula-segunda-derivada-covariante-seccion}.

For the third order, regard $\nabla^2\mathbf{u}$ as a $2$-covariant tensor with values in $\mathbf{E}$. The definition of the induced connection gives
\begin{equation}
\label{eq:tercera-derivada-antes-de-expandir-cap2}
(\nabla^3\mathbf{u})(\mathbf{X},\mathbf{Y},\mathbf{Z})=\nabla^{\mathbf{E}}_{\mathbf{X}}\bigl((\nabla^2\mathbf{u})(\mathbf{Y},\mathbf{Z})\bigr)
-(\nabla^2\mathbf{u})(\nabla^M_{\mathbf{X}}\mathbf{Y},\mathbf{Z})-(\nabla^2\mathbf{u})(\mathbf{Y},\nabla^M_{\mathbf{X}}\mathbf{Z}).
\end{equation}
Expand the three terms separately. By
\eqref{eq:formula-segunda-derivada-covariante-seccion},
\[
\nabla^{\mathbf{E}}_{\mathbf{X}}\bigl((\nabla^2\mathbf{u})(\mathbf{Y},\mathbf{Z})\bigr)
=
\nabla^{\mathbf{E}}_{\mathbf{X}}\nabla^{\mathbf{E}}_{\mathbf{Y}}\nabla^{\mathbf{E}}_{\mathbf{Z}}\mathbf{u}
-
\nabla^{\mathbf{E}}_{\mathbf{X}}\nabla^{\mathbf{E}}_{\nabla^M_{\mathbf{Y}}\mathbf{Z}}\mathbf{u}.
\]
Likewise,
\[
(\nabla^2\mathbf{u})(\nabla^M_{\mathbf{X}}\mathbf{Y},\mathbf{Z})
=\nabla^{\mathbf{E}}_{\nabla^M_{\mathbf{X}}\mathbf{Y}}\nabla^{\mathbf{E}}_{\mathbf{Z}}\mathbf{u}
-\nabla^{\mathbf{E}}_{\nabla^M_{\nabla^M_{\mathbf{X}}\mathbf{Y}}\mathbf{Z}}\mathbf{u},
\]
and
\[
(\nabla^2\mathbf{u})(\mathbf{Y},\nabla^M_{\mathbf{X}}\mathbf{Z})
=\nabla^{\mathbf{E}}_{\mathbf{Y}}\nabla^{\mathbf{E}}_{\nabla^M_{\mathbf{X}}\mathbf{Z}}\mathbf{u}
-\nabla^{\mathbf{E}}_{\nabla^M_{\mathbf{Y}}\nabla^M_{\mathbf{X}}\mathbf{Z}}\mathbf{u}.
\]
Substituting these three identities into
\eqref{eq:tercera-derivada-antes-de-expandir-cap2}, keeping the two negative signs in its second line, yields exactly
\eqref{eq:formula-tercera-derivada-covariante-seccion}.
\end{proof}

\begin{remark}[On the use of normal coordinates at third order]
\label{obs:tercera-derivada-no-se-obtiene-derivando-identidad-puntual}
Even if normal coordinates centered at a point $p$ are chosen, the identity
\[
(\nabla^2\mathbf{u})(\mathbf{X},\mathbf{Y})(p)=\nabla^{\mathbf{E}}_{\mathbf{X}}\nabla^{\mathbf{E}}_{\mathbf{Y}}\mathbf{u}(p)
\]
holds only for extensions of $\mathbf{X}$ and $\mathbf{Y}$ whose first covariant derivatives vanish at $p$. It cannot be differentiated as though it were an equality valid in a neighborhood of $p$. Indeed, formula
\eqref{eq:formula-tercera-derivada-covariante-seccion} contains terms with second covariant derivatives of the auxiliary fields, such as
$\nabla^M_{\mathbf{Y}}\nabla^M_{\mathbf{X}}\mathbf{Z}$, which do not vanish merely because normal coordinates were chosen. For this reason, when contracting or commuting higher-order covariant derivatives, we will first establish tensor identities and only then evaluate them in a normal frame. The operators obtained by such contractions will be defined in the chapter on differential operators on bundles.
\end{remark}

\begin{lemma}[Leibniz rule for interior multiplication]\label{leibniz-multiplicacion-interior-E}\index{Leibniz rule for interior multiplication}
Let $(M,\mathbf{g})$ be a Riemannian manifold with or without boundary equipped with a connection $\nabla$ on $TM$.
Let $(\mathbf{E}_{1},\dots,\mathbf{E}_{n})$ be a smooth local frame of $TM$ with dual coframe
$(\boldsymbol{\varepsilon}^{1},\dots,\boldsymbol{\varepsilon}^{n})$.
Let $\mathbf{E}\to M$ be a smooth vector bundle with connection $\nabla^{\mathbf{E}}$.

For every tensor
\[
\mathbf{F}\in \Gamma\!\bigl(T^{(k,l+1)}(TM)\otimes \mathbf{E}\bigr),
\]
we have
\[
\nabla_{a}\bigl(\iota_{\mathbf{E}_{a}}\mathbf{F}\bigr)
=\iota_{\nabla_{a}\mathbf{E}_{a}}\mathbf{F}+\iota_{\mathbf{E}_{a}}(\nabla_{a}\mathbf{F}),
\]
where $\nabla_{a}=\nabla_{\mathbf{E}_{a}}$ is the total connection on
$T^{(k,l+1)}(TM)\otimes \mathbf{E}$,
and $\iota_{\mathbf{E}_{a}}$ acts only on the covariant block of the tensor.
\end{lemma}

\begin{proof}
The identity holds for each index $a$ before Einstein summation is performed. Fix a point $p$, and evaluate both sides on arbitrary covectors and vectors in the remaining indices of $\mathbf F$. Extend these arguments locally so that their covariant derivatives in the direction $\mathbf E_a$ vanish at $p$. By the definitions of the induced connection and contraction,
\[
\bigl(\iota_{\mathbf E_a}\mathbf F\bigr)
(\boldsymbol\omega^1,\dots,\boldsymbol\omega^k,
 \mathbf X_1,\dots,\mathbf X_l)
=
\mathbf F(\boldsymbol\omega^1,\dots,\boldsymbol\omega^k,
 \mathbf E_a,\mathbf X_1,\dots,\mathbf X_l).
\]
Differentiating this equality in the direction $\mathbf E_a$ and evaluating at $p$, the Leibniz rule for the connection produces exactly two terms:
\[
\begin{split}
\nabla_a\bigl(\iota_{\mathbf E_a}\mathbf F\bigr)
(\boldsymbol\omega^1,\dots,\mathbf X_l)
={}&
\mathbf F(\boldsymbol\omega^1,\dots,
\nabla_a\mathbf E_a,\mathbf X_1,\dots,\mathbf X_l)\\
&+
(\nabla_a\mathbf F)(\boldsymbol\omega^1,\dots,
\mathbf E_a,\mathbf X_1,\dots,\mathbf X_l).
\end{split}
\]
The first summand is
$\bigl(\iota_{\nabla_a\mathbf E_a}\mathbf F\bigr)
(\boldsymbol\omega^1,\dots,\mathbf X_l)$ and the second is
$\bigl(\iota_{\mathbf E_a}\nabla_a\mathbf F\bigr)
(\boldsymbol\omega^1,\dots,\mathbf X_l)$. Since the point and arguments were arbitrary, we obtain the tensor identity; summing over $a$ gives the formula in the statement.
\end{proof}

 Consider the following Leibniz rule:
 \begin{proposition}[Leibniz rule for $(0,2)$ tensors]\label{leibniz tensores (0,2)}\index{Leibniz rule for tensors}
Let $M$ be a smooth manifold with or without boundary, and let
$\mathbf{E}\to M$ be a smooth vector bundle with connection $\nabla^{\mathbf{E}}$.
If $\mathbf{T}\in\Gamma(\mathbf{E}^*\otimes \mathbf{E}^*)$, then, for every $\mathbf{u},\mathbf{v}\in\Gamma(\mathbf{E})$ and every $\mathbf{X}\in\mathfrak X(M)$, we have
\[
\mathbf{X}\big(\mathbf{T}(\mathbf{u},\mathbf{v})\big)=(\nabla^{\mathbf{E}^*\otimes \mathbf{E}^*}_{\mathbf{X}} \mathbf{T})(\mathbf{u},\mathbf{v})+\mathbf{T}(\nabla^{\mathbf{E}}_{\mathbf{X}} \mathbf{u},\mathbf{v})+\mathbf{T}(\mathbf{u},\nabla^{\mathbf{E}}_{\mathbf{X}} \mathbf{v}).
\]
\end{proposition}

\begin{proof}
By definition, for $\mathbf{u},\mathbf{v}\in\Gamma(\mathbf{E})$ the function $\mathbf{T}(\mathbf{u},\mathbf{v})\in C^\infty(M)$ is given by
\[
\mathbf{T}(\mathbf{u},\mathbf{v})(p)=\mathbf{T}_p(\mathbf{u}_p,\mathbf{v}_p),\qquad p\in M.
\]

Since $\Gamma(\mathbf{E}^*\otimes \mathbf{E}^*)$ is a $C^\infty(M)$-module generated by simple tensors, it suffices to prove the identity when $\mathbf{T}=\boldsymbol{\alpha}\otimes\boldsymbol{\beta}$ with $\boldsymbol{\alpha},\boldsymbol{\beta}\in\Gamma(\mathbf{E}^*)$ and then conclude by linearity. In this case,
\[
\mathbf{T}(\mathbf{u},\mathbf{v})=\langle \boldsymbol{\alpha},\mathbf{u}\rangle\langle \boldsymbol{\beta},\mathbf{v}\rangle.
\]
Applying the product rule for functions and the characterization of the induced connection on the dual (Proposition~\ref{conexion en tensores, duales y endomorfismos}),
\[
\mathbf{X}\big(\mathbf{T}(\mathbf{u},\mathbf{v})\big)
=\mathbf{X}\big(\langle \boldsymbol{\alpha},\mathbf{u}\rangle\langle \boldsymbol{\beta},\mathbf{v}\rangle\big)
\]
\[
=\mathbf{X}\langle \boldsymbol{\alpha},\mathbf{u}\rangle\langle \boldsymbol{\beta},\mathbf{v}\rangle + \langle \boldsymbol{\alpha},\mathbf{u}\rangle \mathbf{X} \langle \boldsymbol{\beta},\mathbf{v}\rangle
\]
\[
=\big(\langle \nabla^{\mathbf{E}^*}_{\mathbf{X}}\boldsymbol{\alpha},\mathbf{u}\rangle+\langle \boldsymbol{\alpha},\nabla^{\mathbf{E}}_{\mathbf{X}} \mathbf{u}\rangle\big)\langle \boldsymbol{\beta},\mathbf{v}\rangle
+\langle \boldsymbol{\alpha},\mathbf{u}\rangle\big(\langle \nabla^{\mathbf{E}^*}_{\mathbf{X}}\boldsymbol{\beta},\mathbf{v}\rangle+\langle \boldsymbol{\beta},\nabla^{\mathbf{E}}_{\mathbf{X}} \mathbf{v}\rangle\big)
\]
\[
=\underbrace{\langle \nabla^{\mathbf{E}^*}_{\mathbf{X}}\boldsymbol{\alpha},\mathbf{u}\rangle\langle \boldsymbol{\beta},\mathbf{v}\rangle
+\langle \boldsymbol{\alpha},\mathbf{u}\rangle\langle \nabla^{\mathbf{E}^*}_{\mathbf{X}}\boldsymbol{\beta},\mathbf{v}\rangle}_{(\nabla^{\mathbf{E}^*\otimes \mathbf{E}^*}_{\mathbf{X}}(\boldsymbol{\alpha}\otimes\boldsymbol{\beta}))(\mathbf{u},\mathbf{v})\ \text{by Proposition~\ref{conexion en tensores, duales y endomorfismos}.}}
+\underbrace{\langle \boldsymbol{\alpha},\nabla^{\mathbf{E}}_{\mathbf{X}} \mathbf{u}\rangle\langle \boldsymbol{\beta},\mathbf{v}\rangle}_{\mathbf{T}(\nabla^{\mathbf{E}}_{\mathbf{X}} \mathbf{u},\mathbf{v})}
+\underbrace{\langle \boldsymbol{\alpha},\mathbf{u}\rangle\langle \boldsymbol{\beta},\nabla^{\mathbf{E}}_{\mathbf{X}} \mathbf{v}\rangle}_{\mathbf{T}(\mathbf{u},\nabla^{\mathbf{E}}_{\mathbf{X}} \mathbf{v})}.
\]
Consequently,
\[
\mathbf{X}\big(\mathbf{T}(\mathbf{u},\mathbf{v})\big)=(\nabla^{\mathbf{E}^*\otimes \mathbf{E}^*}_{\mathbf{X}} \mathbf{T})(\mathbf{u},\mathbf{v})+\mathbf{T}(\nabla^{\mathbf{E}}_{\mathbf{X}} \mathbf{u},\mathbf{v})+\mathbf{T}(\mathbf{u},\nabla^{\mathbf{E}}_{\mathbf{X}} \mathbf{v})
\]
for $\mathbf{T}=\boldsymbol{\alpha}\otimes\boldsymbol{\beta}$. By linearity of every term, the equality extends to all $\mathbf{T}\in\Gamma(\mathbf{E}^*\otimes \mathbf{E}^*)$.
\end{proof}
This Leibniz rule motivates the following more general concept, which describes when a bundle homomorphism respects the structures given by the connection.
 \begin{definition}[Parallel smooth bundle homomorphisms]\label{def: homomorfismos paralelos}\index{parallel smooth bundle homomorphisms}
Let $M$ be a smooth manifold with or without boundary, and let
$\mathbf{E},\mathbf{F},\mathbf{G} \to M$ be smooth vector bundles with connections
$\nabla^{\mathbf{E}}, \nabla^{\mathbf{F}}, \nabla^{\mathbf{G}}$.
A smooth bundle homomorphism $\Phi\colon \mathbf{E} \otimes \mathbf{F} \longrightarrow \mathbf{G}$ is said to be parallel if, for all $\mathbf{u}\in \Gamma(\mathbf{E})$, $\mathbf{v}\in \Gamma(\mathbf{F})$, and $\mathbf{X}\in \mathfrak{X}(M)$, we have
\[
\nabla^{\mathbf{G}}_{\mathbf{X}}\big(\Phi(\mathbf{u}\otimes \mathbf{v})\big)
= \Phi\big((\nabla^{\mathbf{E}}_{\mathbf{X}} \mathbf{u})\otimes \mathbf{v}\big)
+ \Phi\big(\mathbf{u}\otimes (\nabla^{\mathbf{F}}_{\mathbf{X}} \mathbf{v})\big).
\]
\end{definition}
This can be written equivalently as shown in the following proposition:
\begin{proposition}\label{prop:paralelo-hom}
Let $M$ be a smooth manifold with or without boundary, and let
$\mathbf{E},\mathbf{F},\mathbf{G}\to M$ be smooth vector bundles with connections
$\nabla^{\mathbf{E}},\nabla^{\mathbf{F}},\nabla^{\mathbf{G}}$. Let
\[
\Phi\colon \mathbf{E}\otimes \mathbf{F}\longrightarrow \mathbf{G}
\]
be a smooth bundle homomorphism. Then $\Phi$ is a parallel smooth bundle homomorphism in the sense of Definition~\ref{def: homomorfismos paralelos} if and only if, when regarded as a section of $\operatorname{Hom}(\mathbf{E}\otimes \mathbf{F},\mathbf{G})$, it satisfies $\nabla_{\mathbf{X}} \Phi=0$,
$\forall \mathbf{X}\in \mathfrak{X}(M)$, where $\nabla$ is the induced connection on $\operatorname{Hom}(\mathbf{E}\otimes \mathbf{F},\mathbf{G})$.
\end{proposition}

\begin{proof}
By the definition of the connection on $\operatorname{Hom}(\mathbf{E}\otimes \mathbf{F},\mathbf{G})$, we have
\[
(\nabla_{\mathbf{X}}\Phi)(\mathbf{u}\otimes \mathbf{v})
=
\nabla^{\mathbf{G}}_{\mathbf{X}}\bigl(\Phi(\mathbf{u}\otimes \mathbf{v})\bigr)
-
\Phi\bigl(\nabla_{\mathbf{X}}(\mathbf{u}\otimes \mathbf{v})\bigr),
\]
for all $\mathbf{u}\in\Gamma(\mathbf{E})$, $\mathbf{v}\in\Gamma(\mathbf{F})$, and $\mathbf{X}\in\mathfrak{X}(M)$.

Since $\nabla_{\mathbf{X}}(\mathbf{u}\otimes \mathbf{v})=(\nabla^{\mathbf{E}}_{\mathbf{X}} \mathbf{u})\otimes \mathbf{v} + \mathbf{u}\otimes (\nabla^{\mathbf{F}}_{\mathbf{X}} \mathbf{v})$, we obtain
\[
(\nabla_{\mathbf{X}}\Phi)(\mathbf{u}\otimes \mathbf{v})
=
\nabla^{\mathbf{G}}_{\mathbf{X}}\bigl(\Phi(\mathbf{u}\otimes \mathbf{v})\bigr)
-
\Phi\big((\nabla^{\mathbf{E}}_{\mathbf{X}} \mathbf{u})\otimes \mathbf{v}\big)
-
\Phi\big(\mathbf{u}\otimes (\nabla^{\mathbf{F}}_{\mathbf{X}} \mathbf{v})\big).
\]

Therefore $(\nabla_{\mathbf{X}}\Phi)(\mathbf{u}\otimes \mathbf{v})=0$ for every $\mathbf{u},\mathbf{v}$ if and only if
\[
\nabla^{\mathbf{G}}_{\mathbf{X}}\bigl(\Phi(\mathbf{u}\otimes \mathbf{v})\bigr)
=
\Phi\big((\nabla^{\mathbf{E}}_{\mathbf{X}} \mathbf{u})\otimes \mathbf{v}\big)
+
\Phi\big(\mathbf{u}\otimes (\nabla^{\mathbf{F}}_{\mathbf{X}} \mathbf{v})\big),
\]
that is, if and only if $\Phi$ is a parallel smooth bundle homomorphism in the sense of Definition~\ref{def: homomorfismos paralelos}.

Since these simple tensor fields generate $\Gamma(\mathbf{E}\otimes \mathbf{F})$, the equivalence holds for every section, completing the proof.
\end{proof}
\begin{example}[Evaluation]\label{ej:variedades-riemannianas-la-evaluacion}
Let $M$ be a smooth manifold with or without boundary, let $\mathbf{E} \to M$ be a smooth vector bundle with connection $\nabla^{\mathbf{E}}$, and let
$\nabla^{\mathbf{E}^*}$ be the induced connection on the dual.

Define the smooth bundle homomorphism called the \textit{evaluation homomorphism}
\[
\mathrm{ev}\colon \mathbf{E}^*\otimes \mathbf{E} \longrightarrow M\times \mathbb{R}
\]
pointwise by
\[
\mathrm{ev}_p(\mathbf{v}_p\otimes \mathbf{u}_p)=\langle \mathbf{v}_p,\mathbf{u}_p\rangle,
\qquad p\in M, \mathbf{v}_p\in \mathbf{E}^*_p, \mathbf{u}_p\in \mathbf{E}_p.
\]

The trivial bundle $M\times \mathbb{R}$ carries the trivial connection;
that is, regarding $f\in C^\infty(M)$ as a section of $M\times \mathbb{R}$, we have
\[
\nabla^{M\times \mathbb{R}}_{\mathbf{X}} f = \mathbf{X}(f).
\]
Let $\mathbf{v}\in \Gamma(\mathbf{E}^*)$, $\mathbf{u}\in \Gamma(\mathbf{E})$, and $\mathbf{X}\in \mathfrak{X}(M)$.
On the left-hand side:
\[
\nabla^{M\times \mathbb{R}}_{\mathbf{X}}(\mathrm{ev}(\mathbf{v}\otimes \mathbf{u}))
= \nabla^{M\times \mathbb{R}}_{\mathbf{X}}\big(p\mapsto \langle \mathbf{v}_p,\mathbf{u}_p\rangle\big)
= \mathbf{X}(\langle \mathbf{v},\mathbf{u}\rangle).
\]
On the right-hand side:
\[
\mathrm{ev}\big((\nabla^{\mathbf{E}^*}_{\mathbf{X}} \mathbf{v})\otimes \mathbf{u} + \mathbf{v}\otimes \nabla^{\mathbf{E}}_{\mathbf{X}} \mathbf{u}\big)
= \langle \nabla^{\mathbf{E}^*}_{\mathbf{X}} \mathbf{v},\mathbf{u}\rangle + \langle \mathbf{v},\nabla^{\mathbf{E}}_{\mathbf{X}} \mathbf{u}\rangle.
\]

By the definition of the induced connection on $\mathbf{E}^*$,
\[
\mathbf{X}(\langle \mathbf{v},\mathbf{u}\rangle) = \langle \nabla^{\mathbf{E}^*}_{\mathbf{X}} \mathbf{v},\mathbf{u}\rangle + \langle \mathbf{v},\nabla^{\mathbf{E}}_{\mathbf{X}} \mathbf{u}\rangle.
\]
Thus,
\[
\nabla^{M\times \mathbb{R}}_{\mathbf{X}}\big(\mathrm{ev}(\mathbf{v}\otimes \mathbf{u})\big)
= \mathrm{ev}\big((\nabla^{\mathbf{E}^*}_{\mathbf{X}} \mathbf{v})\otimes \mathbf{u}+ \mathbf{v}\otimes \nabla^{\mathbf{E}}_{\mathbf{X}} \mathbf{u}\big),
\]
and we conclude that evaluation is a parallel smooth bundle homomorphism.
\end{example}

\begin{example}[Contraction]\label{ej:variedades-riemannianas-la-contraccion}
Let $M$ be a smooth manifold with or without boundary, let $\mathbf{E} \to M$ be a smooth vector bundle with connection $\nabla^{\mathbf{E}}$, and let $\mathbf{F}\to M$ be another bundle with connection $\nabla^{\mathbf{F}}$.

Define the smooth homomorphism called the \textit{contraction homomorphism}
\[
\mathrm{contr}\colon \mathbf{E}^*\otimes \mathbf{E} \otimes \mathbf{F} \longrightarrow \mathbf{F}
\]
pointwise by
\[
\mathrm{contr}_p(\mathbf{v}_p\otimes \mathbf{u}_p\otimes \mathbf{w}_p)=\langle \mathbf{v}_p,\mathbf{u}_p\rangle \mathbf{w}_p,
\qquad p\in M, \mathbf{v}_p\in \mathbf{E}^*_p, \mathbf{u}_p\in \mathbf{E}_p, \mathbf{w}_p\in \mathbf{F}_p.
\]

The notion of a parallel homomorphism extends naturally to three factors: a homomorphism
\[
\Phi\colon \mathbf{E}^*\otimes \mathbf{E} \otimes \mathbf{F} \longrightarrow \mathbf{G}
\]
is parallel if, for all $\mathbf{v}\in\Gamma(\mathbf{E}^*)$, $\mathbf{u}\in\Gamma(\mathbf{E})$, $\mathbf{w}\in\Gamma(\mathbf{F})$, and $\mathbf{X}\in\mathfrak{X}(M)$, we have
\[
\nabla^{\mathbf{G}}_{\mathbf{X}}(\Phi(\mathbf{v}\otimes \mathbf{u}\otimes \mathbf{w}))
= \Phi\big((\nabla^{\mathbf{E}^*}_{\mathbf{X}} \mathbf{v})\otimes \mathbf{u}\otimes \mathbf{w}\big)
+ \Phi\big(\mathbf{v}\otimes (\nabla^{\mathbf{E}}_{\mathbf{X}} \mathbf{u})\otimes \mathbf{w}\big)
+ \Phi\big(\mathbf{v}\otimes \mathbf{u}\otimes \nabla^{\mathbf{F}}_{\mathbf{X}} \mathbf{w}\big).
\]
Let us verify that contraction satisfies this property.

For the left-hand side, we have
\[
\begin{aligned}
\nabla^{\mathbf{F}}_{\mathbf{X}}(\mathrm{contr}(\mathbf{v}\otimes \mathbf{u}\otimes \mathbf{w}))
&=\nabla^{\mathbf{F}}_{\mathbf{X}}\big(p\mapsto \langle \mathbf{v}_p,\mathbf{u}_p\rangle \mathbf{w}_p\big) \\[0.5em]
&=\nabla^{\mathbf{F}}_{\mathbf{X}}(\langle \mathbf{v},\mathbf{u}\rangle \cdot \mathbf{w}).
\end{aligned}
\]
The Leibniz rule on $\mathbf{F}$ gives
\[
\nabla^{\mathbf{F}}_{\mathbf{X}}(\langle \mathbf{v},\mathbf{u}\rangle \cdot \mathbf{w})
= \mathbf{X}(\langle \mathbf{v},\mathbf{u}\rangle)\mathbf{w} + \langle \mathbf{v},\mathbf{u}\rangle\nabla^{\mathbf{F}}_{\mathbf{X}} \mathbf{w},
\]
so
\[
\nabla^{\mathbf{F}}_{\mathbf{X}}(\mathrm{contr}(\mathbf{v}\otimes \mathbf{u}\otimes \mathbf{w}))
= \mathbf{X}(\langle \mathbf{v},\mathbf{u}\rangle)\mathbf{w} + \langle \mathbf{v},\mathbf{u}\rangle\nabla^{\mathbf{F}}_{\mathbf{X}} \mathbf{w}.
\]

We wish to calculate
\[
\ \mathrm{contr}\big((\nabla^{\mathbf{E}^*}_{\mathbf{X}} \mathbf{v})\otimes \mathbf{u} \otimes \mathbf{w}\big)
+ \mathrm{contr}\big(\mathbf{v}\otimes (\nabla^{\mathbf{E}}_{\mathbf{X}} \mathbf{u})\otimes \mathbf{w}\big)
+ \mathrm{contr}\big(\mathbf{v}\otimes \mathbf{u} \otimes \nabla^{\mathbf{F}}_{\mathbf{X}} \mathbf{w}\big).
\]
Evaluating each term using the definition of contraction gives
\[
\begin{aligned}
\mathrm{contr}\big((\nabla^{\mathbf{E}^*}_{\mathbf{X}} \mathbf{v})\otimes \mathbf{u} \otimes \mathbf{w}\big) &= \langle \nabla^{\mathbf{E}^*}_{\mathbf{X}} \mathbf{v}, \mathbf{u}\rangle \mathbf{w}, \\[0.5em]
\mathrm{contr}\big(\mathbf{v}\otimes (\nabla^{\mathbf{E}}_{\mathbf{X}} \mathbf{u})\otimes \mathbf{w}\big) &= \langle \mathbf{v}, \nabla^{\mathbf{E}}_{\mathbf{X}} \mathbf{u}\rangle \mathbf{w}, \\[0.5em]
\mathrm{contr}\big(\mathbf{v}\otimes \mathbf{u} \otimes \nabla^{\mathbf{F}}_{\mathbf{X}} \mathbf{w}\big) &= \langle \mathbf{v},\mathbf{u}\rangle\nabla^{\mathbf{F}}_{\mathbf{X}} \mathbf{w}.
\end{aligned}
\]
Thus,
\[
\mathrm{contr}\big((\nabla^{\mathbf{E}^*}_{\mathbf{X}} \mathbf{v})\otimes \mathbf{u} \otimes \mathbf{w}\big)
+ \mathrm{contr}\big(\mathbf{v}\otimes (\nabla^{\mathbf{E}}_{\mathbf{X}} \mathbf{u})\otimes \mathbf{w}\big)
+ \mathrm{contr}\big(\mathbf{v}\otimes \mathbf{u} \otimes \nabla^{\mathbf{F}}_{\mathbf{X}} \mathbf{w}\big)
= \langle \nabla^{\mathbf{E}^*}_{\mathbf{X}} \mathbf{v}, \mathbf{u}\rangle \mathbf{w}
+ \langle \mathbf{v}, \nabla^{\mathbf{E}}_{\mathbf{X}} \mathbf{u}\rangle \mathbf{w}
+ \langle \mathbf{v},\mathbf{u}\rangle\nabla^{\mathbf{F}}_{\mathbf{X}} \mathbf{w}.
\]

Finally, by the definition of the connection on the dual $\mathbf{E}^*$, we have
\[
\mathbf{X}(\langle \mathbf{v},\mathbf{u}\rangle)=\langle \nabla^{\mathbf{E}^*}_{\mathbf{X}} \mathbf{v}, \mathbf{u}\rangle+\langle \mathbf{v},\nabla^{\mathbf{E}}_{\mathbf{X}} \mathbf{u}\rangle.
\]
Equating the two terms, we obtain
\[
\nabla^{\mathbf{F}}_{\mathbf{X}}(\mathrm{contr}(\mathbf{v}\otimes \mathbf{u}\otimes \mathbf{w}))
= \mathrm{contr}\big((\nabla^{\mathbf{E}^*}_{\mathbf{X}} \mathbf{v})\otimes \mathbf{u}\otimes \mathbf{w}\big)
+ \mathrm{contr}\big(\mathbf{v}\otimes (\nabla^{\mathbf{E}}_{\mathbf{X}} \mathbf{u})\otimes \mathbf{w}\big)
+ \mathrm{contr}\big(\mathbf{v}\otimes \mathbf{u} \otimes \nabla^{\mathbf{F}}_{\mathbf{X}} \mathbf{w}\big).
\]
Therefore contraction is a parallel smooth bundle homomorphism.
\end{example}

\begin{proposition}[Compatibility of the connection with contractions]
\label{prop:conexion-conmuta-contracciones}
\index{induced connection!compatibility with contractions}
Let $M$ be a smooth manifold with or without boundary, let
$\mathbf{E}_1,\dots,\mathbf{E}_N,\mathbf{F}\to M$ be smooth vector bundles equipped with connections, and let
\[
C\colon \mathbf{E}_1\otimes\cdots\otimes \mathbf{E}_N\longrightarrow \mathbf{F}
\]
be a parallel smooth bundle homomorphism. Then, for every section
$\mathbf{T}\in\Gamma(\mathbf{E}_1\otimes\cdots\otimes \mathbf{E}_N)$ and every
$\mathbf{X}\in\mathfrak X(M)$,
\begin{equation}
\label{eq:conexion-conmuta-homomorfismo-paralelo}
\nabla^{\mathbf{F}}_{\mathbf{X}}\bigl(C(\mathbf{T})\bigr)=C(\nabla_{\mathbf{X}}\mathbf{T}).
\end{equation}
In particular, induced connections commute with evaluation, permutations of factors, and all canonical contractions.

If $(M,\mathbf{g})$ is Riemannian and $\nabla^M$ is the Levi--Civita connection, metric contraction is parallel. For a tensor $S$ with
$l\geq2$ covariant indices and $1\leq a<b\leq l$, let
$\operatorname{tr}_{\mathbf{g}}^{(a,b)}S$ denote contraction of covariant
slots $a$ and $b$ with $\mathbf{g}^{-1}$. Then
\begin{equation}
\label{eq:derivada-conmuta-traza-metrica}
\nabla\bigl(\operatorname{tr}_{\mathbf{g}}^{(a,b)}S\bigr)
=\operatorname{tr}_{\mathbf{g}}^{(a+1,b+1)}(\nabla S).
\end{equation}
The new derivative index remains free, and the same indices inherited
from $S$ are contracted. In particular, for a covariant $2$-tensor,
\[
 \bigl(\nabla\operatorname{tr}_{\mathbf{g}}S\bigr)_c
 =g^{ij}(\nabla S)_{cij}.
\]
The same compatibility holds for any finite composition of metric
contractions, preserving the inherited slots being contracted.
\end{proposition}

\begin{proof}
By Proposition~\ref{prop:paralelo-hom}, the condition that
$C$ be parallel is equivalent to $\nabla C=0$ when $C$ is regarded as a section of the corresponding homomorphism bundle. The Leibniz rule for evaluation gives
\[
\nabla^{\mathbf{F}}_{\mathbf{X}}(C(\mathbf{T}))=(\nabla_{\mathbf{X}}C)(\mathbf{T})+C(\nabla_{\mathbf{X}}\mathbf{T})=C(\nabla_{\mathbf{X}}\mathbf{T}),
\]
which proves \eqref{eq:conexion-conmuta-homomorfismo-paralelo}.

Evaluation and contraction are parallel by the preceding examples.
Permutations of factors are parallel because the product connection differentiates each factor, while the permutation only changes their order. In the Riemannian case, $\nabla^Mg=0$ and also $\nabla^Mg^{-1}=0$; hence raising indices and contracting them with $g^{-1}$ are parallel homomorphisms. Identity \eqref{eq:derivada-conmuta-traza-metrica} is therefore an instance of
\eqref{eq:conexion-conmuta-homomorfismo-paralelo}. The assertion for a finite composition follows by applying the same argument successively.
\end{proof}

Finally, we give an example showing why parallelism of homomorphisms is more general than the Leibniz rule for tensor fields of type $(0,2)$:
\begin{example}[Evaluation and the Leibniz rule]\label{ej:variedades-riemannianas-la-evaluacion-y-la-regla-de-leibniz}
Let $M$ be a smooth manifold with or without boundary, and let $\mathbf{E}\to M$ be a smooth vector bundle with connection $\nabla^{\mathbf{E}}$. Define the map
\[
\widetilde{\mathrm{ev}}\colon (\mathbf{E}^*\otimes \mathbf{E}^*)\otimes(\mathbf{E}\otimes \mathbf{E})\longrightarrow M\times \mathbb{R},
\qquad
\widetilde{\mathrm{ev}}_p(\beta_p\otimes(\mathbf{u}_p\otimes \mathbf{v}_p))=\langle \beta_p,\mathbf{u}_p\otimes \mathbf{v}_p\rangle=\beta_{p}(\mathbf{u}_{p},\mathbf{v}_{p}).
\]
Then, for all $\mathbf{T}\in\Gamma(\mathbf{E}^*\otimes \mathbf{E}^*)$, $\mathbf{u},\mathbf{v}\in\Gamma(\mathbf{E})$, and $\mathbf{X}\in\mathfrak X(M)$, the Leibniz rule
\[
\mathbf{X}\big(\mathbf{T}(\mathbf{u},\mathbf{v})\big) = (\nabla^{\mathbf{E}^*\otimes \mathbf{E}^*}_{\mathbf{X}} \mathbf{T})(\mathbf{u},\mathbf{v})+\mathbf{T}(\nabla^{\mathbf{E}}_{\mathbf{X}} \mathbf{u},\mathbf{v})+\mathbf{T}(\mathbf{u},\nabla^{\mathbf{E}}_{\mathbf{X}} \mathbf{v})
\]
is equivalent to $\widetilde{\mathrm{ev}}$ being a parallel smooth bundle homomorphism, that is,
\[
\nabla^{M\times\mathbb{R}}_{\mathbf{X}}\big(\widetilde{\mathrm{ev}}(\mathbf{T}\otimes(\mathbf{u}\otimes \mathbf{v}))\big)
=
\widetilde{\mathrm{ev}}\big((\nabla^{\mathbf{E}^*\otimes \mathbf{E}^*}_{\mathbf{X}} \mathbf{T})\otimes(\mathbf{u}\otimes \mathbf{v})\big)
+\widetilde{\mathrm{ev}}\big(\mathbf{T}\otimes((\nabla^{\mathbf{E}}_{\mathbf{X}} \mathbf{u})\otimes \mathbf{v})\big)
+\widetilde{\mathrm{ev}}\big(\mathbf{T}\otimes(\mathbf{u}\otimes \nabla^{\mathbf{E}}_{\mathbf{X}} \mathbf{v})\big).
\]
\end{example}

\begin{proof}
By the pointwise definition,
\[
\widetilde{\mathrm{ev}}(\mathbf{T}\otimes(\mathbf{u}\otimes \mathbf{v})) = \langle \mathbf{T},\mathbf{u}\otimes \mathbf{v}\rangle = \mathbf{T}(\mathbf{u},\mathbf{v}),
\qquad
\nabla^{M\times\mathbb{R}}_{\mathbf{X}} f = \mathbf{X}(f)\ \text{for } f\in C^\infty(M).
\]

If $\widetilde{\mathrm{ev}}$ is parallel, then
\[
\mathbf{X}\big(\mathbf{T}(\mathbf{u},\mathbf{v})\big)
= \nabla^{M\times\mathbb{R}}_{\mathbf{X}}\big(\widetilde{\mathrm{ev}}(\mathbf{T}\otimes(\mathbf{u}\otimes \mathbf{v}))\big)
\]
\[
= \widetilde{\mathrm{ev}}\big((\nabla^{\mathbf{E}^*\otimes \mathbf{E}^*}_{\mathbf{X}} \mathbf{T})\otimes(\mathbf{u}\otimes \mathbf{v})\big)
+ \widetilde{\mathrm{ev}}\big(\mathbf{T}\otimes((\nabla^{\mathbf{E}}_{\mathbf{X}} \mathbf{u})\otimes \mathbf{v})\big)
+ \widetilde{\mathrm{ev}}\big(\mathbf{T}\otimes(\mathbf{u}\otimes \nabla^{\mathbf{E}}_{\mathbf{X}} \mathbf{v})\big)
\]
\[
= \langle \nabla^{\mathbf{E}^*\otimes \mathbf{E}^*}_{\mathbf{X}} \mathbf{T},\mathbf{u}\otimes \mathbf{v}\rangle
+ \langle \mathbf{T},(\nabla^{\mathbf{E}}_{\mathbf{X}} \mathbf{u})\otimes \mathbf{v}\rangle
+ \langle \mathbf{T},\mathbf{u}\otimes \nabla^{\mathbf{E}}_{\mathbf{X}} \mathbf{v}\rangle
\]
\[
= (\nabla^{\mathbf{E}^*\otimes \mathbf{E}^*}_{\mathbf{X}} \mathbf{T})(\mathbf{u},\mathbf{v}) + \mathbf{T}(\nabla^{\mathbf{E}}_{\mathbf{X}} \mathbf{u},\mathbf{v}) + \mathbf{T}(\mathbf{u},\nabla^{\mathbf{E}}_{\mathbf{X}} \mathbf{v}).
\]
This recovers the Leibniz rule.

Conversely, if the Leibniz rule holds for every $\mathbf{T},\mathbf{u},\mathbf{v}$, then
\[
\nabla^{M\times\mathbb{R}}_{\mathbf{X}}\big(\widetilde{\mathrm{ev}}(\mathbf{T}\otimes(\mathbf{u}\otimes \mathbf{v}))\big)
= \mathbf{X}\big(\mathbf{T}(\mathbf{u},\mathbf{v})\big)
\]
\[
= (\nabla^{\mathbf{E}^*\otimes \mathbf{E}^*}_{\mathbf{X}} \mathbf{T})(\mathbf{u},\mathbf{v}) + \mathbf{T}(\nabla^{\mathbf{E}}_{\mathbf{X}} \mathbf{u},\mathbf{v}) + \mathbf{T}(\mathbf{u},\nabla^{\mathbf{E}}_{\mathbf{X}} \mathbf{v})
\]
\[
= \langle \nabla^{\mathbf{E}^*\otimes \mathbf{E}^*}_{\mathbf{X}} \mathbf{T},\mathbf{u}\otimes \mathbf{v}\rangle
+ \langle \mathbf{T},(\nabla^{\mathbf{E}}_{\mathbf{X}} \mathbf{u})\otimes \mathbf{v}\rangle
+ \langle \mathbf{T},\mathbf{u}\otimes \nabla^{\mathbf{E}}_{\mathbf{X}} \mathbf{v}\rangle
\]
\[
= \widetilde{\mathrm{ev}}\big((\nabla^{\mathbf{E}^*\otimes \mathbf{E}^*}_{\mathbf{X}} \mathbf{T})\otimes(\mathbf{u}\otimes \mathbf{v})\big)
+\widetilde{\mathrm{ev}}\big(\mathbf{T}\otimes((\nabla^{\mathbf{E}}_{\mathbf{X}} \mathbf{u})\otimes \mathbf{v})\big)
+\widetilde{\mathrm{ev}}\big(\mathbf{T}\otimes(\mathbf{u}\otimes \nabla^{\mathbf{E}}_{\mathbf{X}} \mathbf{v})\big).
\]
This proves that the Leibniz rule is equivalent to $\widetilde{\mathrm{ev}}$ being parallel.
\end{proof}

The condition that a connection $\nabla$ on a smooth vector bundle $\mathbf{E}$ be compatible with a bundle metric $\mathbf{h}$ can be stated compactly as follows:
\begin{definition}\label{conexion compatible con metrica}\index{connection compatible with a bundle metric}\index{connection!compatible with a bundle metric}\glsadd{conexion-compatible}
Let $M$ be a smooth manifold with or without boundary, let
$\mathbf{E}\to M$ be a smooth vector bundle, let $\nabla$ be a connection on
$\mathbf{E}$, and let $\mathbf{h}$ be a bundle metric on $\mathbf{E}$.
We say that $\nabla$ is compatible with $\mathbf{h}$ if
\[
\nabla^{\mathbf{E}^{*}\otimes \mathbf{E}^{*}}\mathbf{h}=0.
\]
In this case, the metric $\mathbf{h}$ is also said to be \textit{parallel} with respect to $\nabla$.
\end{definition}

 This condition is equivalent to the one given in Proposition~\ref{producto interno en tensores compatible con levi civita}:
 \begin{proposition}\label{prop:variedades-riemannianas-haz-vectorial-suave-conexion-metrica-fibrada}
 Let $M$ be a smooth manifold with or without boundary, let
 $\mathbf{E}\longrightarrow M$ be a smooth vector bundle, let $\nabla$ be a connection on $\mathbf{E}$, and let $\mathbf{h}$ be a bundle metric on $\mathbf{E}$.
 Then $\nabla$ is compatible with $\mathbf{h}$ if and only if, for any
 $\mathbf{u},\mathbf{v}\in \Gamma(\mathbf{E})$ and any
 $\mathbf{X}\in \mathfrak{X}(M)$,
 \[\mathbf{X}\bigl(\mathbf{h}(\mathbf{u},\mathbf{v})\bigr)=\mathbf{h}(\nabla_{\mathbf{X}}\mathbf{u},\mathbf{v})+\mathbf{h}(\mathbf{u},\nabla_{\mathbf{X}}\mathbf{v}).\]
 \end{proposition}
 \begin{proof}
 By Proposition~\ref{leibniz tensores (0,2)}, the covariant derivative of the section $\mathbf h\in\Gamma(T^{(0,2)}(\mathbf E))$ satisfies
 \[
 (\nabla_{\mathbf X}\mathbf h)(\mathbf u,\mathbf v)
 =\mathbf X\bigl(\mathbf h(\mathbf u,\mathbf v)\bigr)
  -\mathbf h(\nabla_{\mathbf X}\mathbf u,\mathbf v)
  -\mathbf h(\mathbf u,\nabla_{\mathbf X}\mathbf v).
 \]
 Metric compatibility is the identity $\nabla\mathbf h=0$. Thus this identity is exactly equivalent to the formula in the statement for all $\mathbf X$, $\mathbf u$, and $\mathbf v$.
 \end{proof}
\section{Hermitian bundles and conjugation of connections}
\label{sec:haces-hermitianos-conexiones-conjugadas}

The bundle metrics introduced at the beginning of the chapter are real.
We now fix their complex version and the precise relationship between conjugation, duals, and connections. Throughout this section, vector fields on the base are real, and Hermitian products are linear in the first variable.

\begin{definition}[Hermitian bundle metric]
\label{def:metrica-hermitiana-fibrada}
Let $M$ be a smooth manifold with or without boundary, and let
$\mathbf{E}\to M$ be a smooth complex vector bundle. A \emph{Hermitian bundle metric} is a smooth family $\mathbf{h}_{\mathbf{E}}$ of positive definite Hermitian forms on the fibers, linear in the first variable and antilinear in the second.
Equivalently, the bilinear map
\[
 \widetilde{\mathbf{h}}_{\mathbf{E}}\colon \mathbf{E}\times_M\overline{\mathbf{E}}\longrightarrow
 \underline{\mathbb C}_M,
 \qquad
 \widetilde{\mathbf{h}}_{\mathbf{E}}(v,\overline w)=\mathbf{h}_{\mathbf{E}}(v,w),
\]
is a smooth section of $\mathbf{E}^*\otimes(\overline{\mathbf{E}})^*$ and is positive definite in the sense that $\widetilde{\mathbf{h}}_{\mathbf{E}}(v,\overline v)>0$ for $v\ne0$.
\end{definition}

\begin{proposition}[Metric on the conjugate bundle]
\label{prop:metrica-haz-conjugado}
The formula
\begin{equation}
\label{eq:metrica-haz-conjugado}
 \mathbf{h}_{\overline{\mathbf{E}}}(\overline v,\overline w)
 :=\overline{\mathbf{h}_{\mathbf{E}}(v,w)}
\end{equation}
defines a smooth Hermitian metric on $\overline{\mathbf{E}}$, linear in the first variable and antilinear in the second. Under the identification
$\overline{\overline{\mathbf{E}}}\cong \mathbf{E}$, the twice-conjugated metric is $\mathbf{h}_{\mathbf{E}}$.
\end{proposition}

\begin{proof}
If $\lambda,\mu\in\mathbb C$, the scalar structure of $\overline{\mathbf{E}}$ and the convention for $\mathbf{h}_{\mathbf{E}}$ give
\begin{align*}
 \mathbf{h}_{\overline{\mathbf{E}}}(\lambda\overline v,\mu\overline w)
 &=\overline{\mathbf{h}_{\mathbf{E}}(\overline\lambda v,\overline\mu w)}
 =\lambda\overline\mu\,
   \overline{\mathbf{h}_{\mathbf{E}}(v,w)}.
\end{align*}
Hermitian symmetry, positivity, and smoothness follow by conjugating the corresponding properties of $\mathbf{h}_{\mathbf{E}}$. Conjugating twice recovers the original formula.
\end{proof}

\begin{proposition}[Fiberwise Riesz isomorphism]
\label{prop:riesz-fibrado-haz-conjugado}
The map
\begin{equation}
\label{eq:riesz-fibrado-haz-conjugado}
 \mathcal R_{\mathbf{E}}\colon\overline{\mathbf{E}}\longrightarrow \mathbf{E}^*,
 \qquad
 \mathcal R_{\mathbf{E}}(\overline v)(w):=\mathbf{h}_{\mathbf{E}}(w,v),
\end{equation}
is a smooth complex-linear bundle isomorphism compatible with restrictions. For each $\overline v$, the functional
$w\mapsto \mathbf{h}_{\mathbf{E}}(w,v)$ is complex-linear. In a unitary local frame
$(\mathbf{e}_1,\ldots,\mathbf{e}_r)$, if
$\displaystyle \ell=\displaystyle\sum_{a=1}^{r}\ell_a \mathbf{e}^a$, then
\begin{equation}
\label{eq:inversa-riesz-fibrado-unitario}
 \mathcal R_{\mathbf{E}}^{-1}(\ell)
 =\overline{\displaystyle\sum_{a=1}^r\overline{\ell_a}\,\mathbf{e}_a}.
\end{equation}
In particular, the map $\mathbf{E}\to \mathbf{E}^*$, $v\mapsto \mathbf{h}_{\mathbf{E}}(\,\cdot\,,v)$, is antilinear, whereas
$\mathcal R_{\mathbf{E}}\colon\overline{\mathbf{E}}\to \mathbf{E}^*$ is complex-linear. In the real case, conjugation disappears and we recover the usual real identification.
\end{proposition}

\begin{proof}
Linearity in $w$ is part of the Hermitian convention. For linearity in $\overline v$, using
$\lambda\overline v=\overline{\overline\lambda v}$, we obtain
\[
 \mathcal R_{\mathbf{E}}(\lambda\overline v)(w)
 =\mathbf{h}_{\mathbf{E}}(w,\overline\lambda v)
 =\lambda \mathbf{h}_{\mathbf{E}}(w,v).
\]
In a unitary frame,
$\displaystyle \mathbf{h}_{\mathbf{E}}(\displaystyle\sum_{a=1}^{r}w^a\mathbf{e}_a,\displaystyle\sum_{a=1}^{r}v^a\mathbf{e}_a)
=\displaystyle\sum_{a=1}^{r}w^a\overline{v^a}$; this directly gives \eqref{eq:inversa-riesz-fibrado-unitario}. The formula proves both bijectivity and local smoothness of $\mathcal R_{\mathbf{E}}$ and its inverse. Since the definition is intrinsic, the local formulas glue and commute with every restriction to an open set. The assertion about
$\mathbf{E}\to \mathbf{E}^*$ follows from
$\mathbf{h}_{\mathbf{E}}(w,\lambda v)=\overline\lambda \mathbf{h}_{\mathbf{E}}(w,v)$.
\end{proof}

\begin{definition}[Complex, dual, and conjugate connections]
\label{def:conexion-compleja-dual-conjugada}
Let $M$ be a smooth manifold with or without boundary, and let
$\mathbf{E}\longrightarrow M$ be a smooth complex vector bundle. A
\emph{complex connection} on $\mathbf{E}$ is a map
\[
 \nabla^{\mathbf{E}}\colon\mathfrak X(M)\times\Gamma(\mathbf{E})\longrightarrow\Gamma(\mathbf{E})
\]
that is $C^\infty(M,\mathbb R)$-linear in the real vector field, additive in the section, satisfies
\[
 \nabla^{\mathbf{E}}_{\mathbf{X}}(f\mathbf{u})
 =\mathbf{X}(f)\mathbf{u}+f\nabla^{\mathbf{E}}_{\mathbf{X}}\mathbf{u},
 \qquad f\in C^\infty(M,\mathbb C),
\]
and is complex-linear in $\mathbf{u}$ for fixed $\mathbf{X}$.

The dual connection $\nabla^{\mathbf{E}^*}$ and the conjugate connection
$\nabla^{\overline{\mathbf{E}}}$ are defined, respectively, by
\begin{align}
 (\nabla_{\mathbf{X}}^{\mathbf{E}^*}\boldsymbol{\ell})(\mathbf{u})
 &:=\mathbf{X}(\boldsymbol{\ell}(\mathbf{u}))
   -\boldsymbol{\ell}(\nabla_{\mathbf{X}}^{\mathbf{E}}\mathbf{u}),
 \label{eq:conexion-dual-compleja}\\
 \nabla_{\mathbf{X}}^{\overline{\mathbf{E}}}\overline{\mathbf{u}}
 &:=\overline{\nabla_{\mathbf{X}}^{\mathbf{E}}\mathbf{u}}.
 \label{eq:conexion-conjugada}
\end{align}
\end{definition}

\begin{proposition}[Well-definedness and naturality]
\label{prop:conexion-conjugada-dual-bien-definida}
Let $M$ be a smooth manifold with or without boundary, and let
$\mathbf{E},\mathbf{F}\longrightarrow M$ be smooth complex vector bundles equipped with complex connections. Formulas
\eqref{eq:conexion-dual-compleja} and
\eqref{eq:conexion-conjugada} define complex connections. They commute with restrictions, and double conjugation recovers $\nabla^{\mathbf{E}}$. If
$\mathbf{A}\colon \mathbf{E}\to \mathbf{F}$ is parallel, that is,
$\nabla^{\mathbf{F}}_{\mathbf{X}}(\mathbf{A}\mathbf{u})
=\mathbf{A}(\nabla^{\mathbf{E}}_{\mathbf{X}}\mathbf{u})$, then $\overline{\mathbf{A}}$ is parallel for the conjugate connections. In conjugate frames, the matrix of
$\nabla^{\overline{\mathbf{E}}}$ is the conjugate of the matrix of $\nabla^{\mathbf{E}}$.
\end{proposition}

\begin{proof}
Only the Leibniz rule for the conjugate connection requires attention. If
$f\in C^\infty(M,\mathbb C)$, then
$f\overline{\mathbf{u}}=\overline{\overline f\,\mathbf{u}}$, and, since
$\mathbf{X}$ is real,
$\mathbf{X}(\overline f)=\overline{\mathbf{X}(f)}$. Hence
\begin{align*}
 \nabla_{\mathbf{X}}^{\overline{\mathbf{E}}}(f\overline{\mathbf{u}})
 &=\overline{\nabla_{\mathbf{X}}^{\mathbf{E}}(\overline f\,\mathbf{u})}
 =\overline{\overline{\mathbf{X}(f)}\mathbf{u}
 +\overline f\,\nabla_{\mathbf{X}}^{\mathbf{E}}\mathbf{u}}\\
 &=\mathbf{X}(f)\overline{\mathbf{u}}
 +f\nabla_{\mathbf{X}}^{\overline{\mathbf{E}}}\overline{\mathbf{u}}.
\end{align*}
The properties of the dual connection follow by expanding
$\mathbf{X}(\boldsymbol{\ell}(f\mathbf{u}))$ and using the Leibniz rule for $\nabla^{\mathbf{E}}$; canceling
$\mathbf{X}(f)\boldsymbol{\ell}(\mathbf{u})$ leaves
$(\nabla_{\mathbf{X}}^{\mathbf{E}^*}\boldsymbol{\ell})(f\mathbf{u})
=f(\nabla_{\mathbf{X}}^{\mathbf{E}^*}\boldsymbol{\ell})(\mathbf{u})$.
The remaining assertions are verified directly on local sections.
\end{proof}

\begin{definition}[Hermitian connection]
\label{def:conexion-hermitiana-compatible}
Let $M$ be a smooth manifold with or without boundary, let
$\mathbf{E}\longrightarrow M$ be a smooth complex vector bundle with Hermitian bundle metric $\mathbf{h}_{\mathbf{E}}$, and let $\nabla^{\mathbf{E}}$ be a complex connection. The connection $\nabla^{\mathbf{E}}$ is
\emph{compatible with $\mathbf{h}_{\mathbf{E}}$}, or
\emph{Hermitian}, if, for all $\mathbf{X}\in\mathfrak X(M)$ and
$\mathbf{u},\mathbf{v}\in\Gamma(\mathbf{E})$,
\begin{equation}
\label{eq:compatibilidad-conexion-hermitiana}
 \mathbf{X}\bigl(\mathbf{h}_{\mathbf{E}}(\mathbf{u},\mathbf{v})\bigr)
 =\mathbf{h}_{\mathbf{E}}(\nabla_{\mathbf{X}}^{\mathbf{E}}\mathbf{u},\mathbf{v})
 +\mathbf{h}_{\mathbf{E}}(\mathbf{u},\nabla_{\mathbf{X}}^{\mathbf{E}}\mathbf{v}).
\end{equation}
\end{definition}

\begin{proposition}[Parallelism of the Riesz isomorphism]
\label{prop:riesz-paralelo-conexion-conjugada}
Let $M$ be a smooth manifold with or without boundary, and let
$\mathbf{E}\longrightarrow M$ be a smooth complex vector bundle with Hermitian bundle metric $\mathbf{h}_{\mathbf{E}}$ and complex connection
$\nabla^{\mathbf{E}}$. The connection $\nabla^{\mathbf{E}}$ is Hermitian if and only if the Riesz isomorphism is parallel:
\begin{equation}
\label{eq:riesz-paralelo}
 \nabla_{\mathbf{X}}^{\mathbf{E}^*}
 \bigl(\mathcal R_{\mathbf{E}}(\overline{\mathbf{v}})\bigr)
 =\mathcal R_{\mathbf{E}}\bigl(
 \nabla_{\mathbf{X}}^{\overline{\mathbf{E}}}\overline{\mathbf{v}}\bigr),
 \qquad
 \mathbf{X}\in\mathfrak X(M),\quad \mathbf{v}\in\Gamma(\mathbf{E}).
\end{equation}
In this case, $\nabla^{\overline{\mathbf{E}}}$ is compatible with
$\mathbf{h}_{\overline{\mathbf{E}}}$, and $\mathcal R_{\mathbf{E}}^{-1}$ is also parallel.
\end{proposition}

\begin{proof}
Evaluating the left-hand side of \eqref{eq:riesz-paralelo} on $\mathbf{u}$ and using the dual connection gives
\[
 \mathbf{X}(\mathbf{h}_{\mathbf{E}}(\mathbf{u},\mathbf{v}))
 -\mathbf{h}_{\mathbf{E}}(\nabla_{\mathbf{X}}^{\mathbf{E}}\mathbf{u},\mathbf{v}).
\]
The right-hand side evaluated on $\mathbf{u}$ is
$\mathbf{h}_{\mathbf{E}}(\mathbf{u},\nabla_{\mathbf{X}}^{\mathbf{E}}\mathbf{v})$.
Consequently, equality for all
$\mathbf{u},\mathbf{v},\mathbf{X}$ is exactly equivalent to
\eqref{eq:compatibilidad-conexion-hermitiana}. Conjugating this identity proves compatibility of the conjugate connection with
$\mathbf{h}_{\overline{\mathbf{E}}}$. Finally, differentiating
$\mathcal R_{\mathbf{E}}^{-1}\mathcal R_{\mathbf{E}}=\operatorname{id}_{\overline{\mathbf{E}}}$ and using
\eqref{eq:riesz-paralelo} proves parallelism of the inverse.
\end{proof}

\section{Principal bundles and principal connections}
\label{sec:haces-principales-conexiones-principales}

A vector bundle describes vectors varying over the base, but a local trivialization also requires choosing a frame in each fiber. The various frames are related by the structure group, and this information is not represented by any particular vector in the bundle. A principal bundle records all admissible frames simultaneously, together with the group action that transforms them.

This description also clarifies why a connection is needed. To compare data in different fibers, one must specify which displacements in the total space are regarded as horizontal. The resulting parallel transport generally depends on the path; curvature measures this dependence infinitesimally. Later, principal connections and their curvature will be the natural variables of Yang--Mills theory, while gauge transformations will express changes of frame.

Throughout this section, all principal bundles carry a right action.
The conventions are chosen so that, if $\mathbf{s}_\beta=\mathbf{s}_\alpha
g_{\alpha\beta}$ are two local sections, the connection potentials satisfy
\[
 \mathbf{A}_\beta
 =
 \operatorname{Ad}_{g_{\alpha\beta}^{-1}}\mathbf{A}_\alpha
 +g_{\alpha\beta}^{-1}dg_{\alpha\beta}.
\]
This is also the convention that will later produce the right gauge action
$\mathbf{A}\mathbin{\triangleleft}\mathbf{u}
=\mathbf{u}^{-1}\mathbf{A}\mathbf{u}+\mathbf{u}^{-1}d\mathbf{u}$.
For a compatible external treatment of principal bundles and their connections, see \cite[Ch.~15]{BleeckerBoossIndex}; the transition, curvature, and gauge identities needed in this book are verified below using the right-action convention fixed above.

\begin{definition}[Principal bundle]
\label{def:haz-principal-derecho}
\index{principal bundle}
Let $G$ be a Lie group and $M$ a smooth manifold with or without boundary. A
\textit{right principal bundle with structure group $G$} consists of a smooth manifold $\mathbf{P}$, a smooth surjective map $\pi\colon \mathbf{P}\to M$, and a smooth right action
\[
 R\colon \mathbf{P}\times G\longrightarrow \mathbf{P},\qquad (p,g)\longmapsto p\cdot g,
\]
such that:

\begin{enumerate}[label=(\alph*)]
\item $\pi(p\cdot g)=\pi(p)$ for every $p\in \mathbf{P}$ and $g\in G$;

\item the action is free;

\item each $x\in M$ has an open neighborhood $U\ni x$ and a
$G$-equivariant diffeomorphism
\[
 \Phi\colon\pi^{-1}(U)\longrightarrow U\times G,
 \qquad
 \Phi(p\cdot g)=\bigl(\pi(p),\operatorname{pr}_2(\Phi(p))g\bigr),
\]
such that $\operatorname{pr}_1\circ\Phi=\pi$.
\end{enumerate}
\end{definition}

Every trivialization determines a local section
$\mathbf{s}\colon U\to \mathbf{P}$, $\mathbf{s}(x)=\Phi^{-1}(x,e)$, and each $p\in\pi^{-1}(U)$ can be written uniquely as
\[
 p=\mathbf{s}(x)\cdot g,\qquad x=\pi(p).
\]
If $\mathbf{s}_\alpha$ and $\mathbf{s}_\beta$ are local sections, there exists a unique smooth map
\[
 g_{\alpha\beta}\colon U_\alpha\cap U_\beta\longrightarrow G
\]
such that
\begin{equation}
\label{eq:transicion-secciones-principales}
 \mathbf{s}_\beta=\mathbf{s}_\alpha g_{\alpha\beta}.
\end{equation}

\begin{proposition}[Cocycles, reconstruction, and triviality]
\label{prop:cociclos-reconstruccion-haces-principales}
Let $M$ be a smooth manifold with or without boundary, and let
$\mathbf{P}\to M$ be a right principal bundle. The transition functions satisfy
\[
 g_{\alpha\alpha}=e,
 \qquad
 g_{\beta\alpha}=g_{\alpha\beta}^{-1},
 \qquad
 g_{\alpha\gamma}=g_{\alpha\beta}g_{\beta\gamma}
\]
on their respective domains. Conversely, a cover
$(U_\alpha)$ and smooth maps satisfying these identities determine a right principal bundle up to principal bundle isomorphism: glue the
$U_\alpha\times G$ by
\begin{equation}
\label{eq:pegado-haz-principal-cociclo}
 (x,h)_\beta\sim(x,g_{\alpha\beta}(x)h)_\alpha.
\end{equation}
Moreover, a principal bundle is trivial if and only if it admits a smooth global section.
\end{proposition}

\begin{proof}
Uniqueness in \eqref{eq:transicion-secciones-principales} gives
$g_{\alpha\alpha}=e$ and $g_{\beta\alpha}=g_{\alpha\beta}^{-1}$.
On a triple overlap,
\[
\mathbf s_\gamma=\mathbf s_\beta g_{\beta\gamma}
=\mathbf s_\alpha g_{\alpha\beta}g_{\beta\gamma}
\]
gives $g_{\alpha\gamma}=g_{\alpha\beta}g_{\beta\gamma}$.

Conversely, take the quotient of
$\coprod_\alpha(U_\alpha\times G)$ by the stated relation.
The three cocycle identities prove reflexivity, symmetry, and
transitivity, respectively. The projection
$\pi([x,h]_\alpha)=x$ is well defined and continuous. For each
$\alpha$, the map
\[
q_\alpha\colon U_\alpha\times G\longrightarrow\pi^{-1}(U_\alpha),
\qquad q_\alpha(x,h)=[x,h]_\alpha,
\]
is bijective: every representative with base point $x\in U_\alpha$
is transformed to the unique representative in chart $\alpha$ by
the cocycle. It is a homeomorphism. Indeed, for an open subset
$B\subseteq U_\alpha\times G$, its saturation in the summand
$U_\beta\times G$ is
\[
\{(x,h)\mid x\in U_\alpha\cap U_\beta,
\ (x,g_{\alpha\beta}(x)h)\in B\},
\]
which is open. The quotient topology therefore makes $q_\alpha$
open; continuity follows from continuity of the quotient map.

The quotient is Hausdorff: points over different base points can
be separated by disjoint open subsets of $M$; two points over the
same base point lie in one open set $\pi^{-1}(U_\alpha)$, which is
homeomorphic to the Hausdorff space $U_\alpha\times G$. It is also
second countable: a countable subcover of $(U_\alpha)$ covers $M$,
and the corresponding products have countable bases whose images
form a global countable basis. On overlaps,
\[
q_\alpha^{-1}\circ q_\beta(x,h)
=(x,g_{\alpha\beta}(x)h)
\]
is a diffeomorphism, with inverse determined by $g_{\beta\alpha}$.
These product charts therefore define a smooth manifold whose
boundary corresponds to the boundary of the base. The action
$[x,h]_\alpha\cdot k=[x,hk]_\alpha$ is well defined, smooth, free,
and transitive on each fiber; the $q_\alpha$ are principal
trivializations. Two realizations with these cocycles are identified
by the map that is the identity in coordinates $(x,h)$: the cocycles
ensure agreement on overlaps. This uniqueness is relative to the
prescribed trivializations, since other gauge automorphisms may exist.

If a global smooth section $\mathbf s$ exists, each $p\in\mathbf P_x$
is uniquely $p=\mathbf s(x)g$. The map
$(x,g)\mapsto\mathbf s(x)g$ is a smooth bijection. In a local
trivialization, if $\mathbf s(x)=(x,b(x))$, its inverse is
$(x,h)\mapsto(x,b(x)^{-1}h)$, which is smooth. Hence it is a global
principal trivialization. Conversely, a trivialization carries
$x\mapsto(x,e)$ to a smooth section of $\mathbf P$.
\end{proof}

\begin{definition}[Principal bundle morphisms and pullbacks]
\label{def:morfismo-pullback-haz-principal}
Let $M$ and $N$ be smooth manifolds with or without boundary, let
$\mathbf{P}\to M$ and $\mathbf{Q}\to N$ be principal bundles with the same group
$G$, and let $f\colon M\to N$ be a smooth map.
A morphism between principal bundles with the same group is a smooth map
$\Phi\colon \mathbf{P}\to \mathbf{Q}$ covering $f$ and satisfying
$\Phi(pg)=\Phi(p)g$. The pullback of $\mathbf{Q}\to N$ by $f$ is
\[
 f^*\mathbf{Q}:=\{(x,q)\in M\times \mathbf{Q}\mid f(x)=\pi_{\mathbf{Q}}(q)\},
 \qquad (x,q)g=(x,qg).
\]
With the projection onto $M$, it is a principal bundle, and
$(x,q)\mapsto q$ is the universal morphism covering $f$.
\end{definition}

\begin{proof}
A trivialization of $\mathbf{Q}$ over $V\subseteq N$ induces
$f^{-1}(V)\times G\cong(f^*\mathbf{Q})|_{f^{-1}(V)}$, proving the local assertion. If $\Psi\colon \mathbf{R}\to \mathbf{Q}$ covers
$f\circ\pi_{\mathbf{R}}$, the unique morphism
$\mathbf{R}\to f^*\mathbf{Q}$ over $M$ projecting to $\Psi$ is
$r\mapsto(\pi_{\mathbf{R}}(r),\Psi(r))$.
\end{proof}

\begin{definition}[Associated bundle]
\label{def:haz-asociado-principal}
\index{associated bundle}
Let $M$ be a smooth manifold with or without boundary, let
$\mathbf{P}\to M$ be a right principal bundle with group $G$, and let
$\rho\colon G\to GL(V)$ be a finite-dimensional real or complex linear representation. Define
\[
 \mathbf{P}\times_\rho V:=(\mathbf{P}\times V)/{\sim},
\]
where
\[
 (p\cdot g,v)\sim(p,\rho(g)v).
\]
We write $[p,v]_\rho$ for the class of $(p,v)$. The projection
\[
 [p,v]_\rho\longmapsto\pi(p)
\]
makes $\mathbf{P}\times_\rho V$ a smooth vector bundle over $M$.
\end{definition}

\begin{proposition}[Associated trivializations and cocycles]
\label{prop:trivializaciones-cociclos-haz-asociado}
Let $M$ be a smooth manifold with or without boundary, and let
$\mathbf{P}\to M$ be a right principal bundle. If $\mathbf{s}_\alpha$ is a local principal section, then
\[
 \Psi_\alpha\colon (\mathbf{P}\times_\rho V)|_{U_\alpha}
 \longrightarrow U_\alpha\times V,
 \qquad
 \Psi_\alpha([\mathbf{s}_\alpha(x)g,v])=(x,\rho(g)v),
\]
is well defined and is a vector bundle trivialization. If
$\mathbf{s}_\beta=\mathbf{s}_\alpha g_{\alpha\beta}$, the coordinates are related by
\[
 v_\beta=\rho(g_{\alpha\beta}^{-1})v_\alpha.
\]
Thus the associated bundle depends only on the principal bundle and the representation, and is independent of the chosen local sections.
\end{proposition}

\begin{proof}
The relation $(pg,v)\sim(p,\rho(g)v)$ shows that
$[\mathbf{s}_\alpha gh,v]=[\mathbf{s}_\alpha g,\rho(h)v]$; both representatives have the same image under $\Psi_\alpha$. The inverse is
$(x,v)\mapsto[\mathbf{s}_\alpha(x),v]$. Finally,
$[\mathbf{s}_\beta,v_\beta]=[\mathbf{s}_\alpha,\rho(g_{\alpha\beta})v_\beta]$, so
$v_\alpha=\rho(g_{\alpha\beta})v_\beta$, which is equivalent to the stated formula.
\end{proof}

\begin{proposition}
\label{prop:secciones-haz-asociado-equivariantes}
Let $M$ be a smooth manifold with or without boundary, and let
$\mathbf{E}=\mathbf{P}\times_\rho V\to M$ be an associated bundle. Smooth sections of $\mathbf{E}$ are in bijective correspondence with smooth maps $\widetilde u\colon \mathbf{P}\to V$ satisfying
\begin{equation}
\label{eq:equivarianza-mapeo-haz-asociado}
 \widetilde u(p\cdot g)=\rho(g^{-1})\widetilde u(p).
\end{equation}
The correspondence is given by
\[
 \mathbf{u}(x)=[p,\widetilde u(p)]_\rho,\qquad p\in \mathbf{P}_x.
\]
\end{proposition}

\begin{proof}
If $p$ is replaced by $p\cdot g$, then
\[
 [p\cdot g,\widetilde u(p\cdot g)]_\rho
 =
 [p,\rho(g)\rho(g^{-1})\widetilde u(p)]_\rho
 =
 [p,\widetilde u(p)]_\rho.
\]
Thus the formula defines a section. Conversely, given $\mathbf{u}(x)$ and
$p\in \mathbf{P}_x$, there is a unique $v\in V$ such that $\mathbf{u}(x)=[p,v]_\rho$; set
$\widetilde u(p)=v$. Uniqueness of $v$ gives
\eqref{eq:equivarianza-mapeo-haz-asociado}. Smoothness in both directions is checked in a principal trivialization.
\end{proof}

\begin{definition}[Adjoint and adjoint group bundles]
\label{def:adP-AdP-distincion}
\index{adjoint bundle}
Let $M$ be a smooth manifold with or without boundary, and let
$\mathbf{P}\to M$ be a right principal bundle with group $G$.
Define two distinct bundles:

\begin{enumerate}[label=(\alph*)]
\item The \textit{adjoint group bundle}
\[
 \operatorname{Ad}(\mathbf{P}):=\mathbf{P}\times_{\operatorname{Conj}}G,
 \qquad
 \operatorname{Conj}_g(h)=ghg^{-1}.
\]
Its fibers are groups; in general, it is not a vector bundle.

\item The \textit{adjoint bundle}
\[
 \operatorname{ad}(\mathbf{P}):=\mathbf{P}\times_{\operatorname{Ad}}\mathfrak g.
\]
This is a vector bundle, of rank $\dim\mathfrak g$.
\end{enumerate}
\end{definition}

Fiberwise multiplication on $\operatorname{Ad}(\mathbf{P})$ is well defined. Its sections $\mathbf{u}$ are identified with smooth maps
$\widehat u\colon \mathbf{P}\to G$ satisfying
\begin{equation}
\label{eq:equivarianza-transformacion-gauge}
 \widehat u(p\cdot g)=g^{-1}\widehat u(p)g.
\end{equation}
By contrast, sections $\boldsymbol{\xi}$ of
$\operatorname{ad}(\mathbf{P})$ are identified with maps
$\widehat\xi\colon \mathbf{P}\to\mathfrak g$ such that
\[
 \widehat\xi(p\cdot g)=\operatorname{Ad}_{g^{-1}}\widehat\xi(p).
\]

\begin{definition}[Horizontal and equivariant forms]
\label{def:formas-horizontales-equivariantes-principal}
Let $M$ be a smooth manifold with or without boundary, let
$\mathbf{P}\to M$ be a right principal bundle with group $G$, and let
$\rho\colon G\to GL(V)$ be a representation. A form
$\boldsymbol{\alpha}\in\Omega^k(\mathbf{P};V)$ is \textit{horizontal} if it vanishes whenever at least one of its arguments belongs to
\[
 \mathbf{V}_p\mathbf{P}:=\ker(d\pi_p).
\]
It is \textit{$\rho$-equivariant} if
\[
 R_g^*\boldsymbol{\alpha}=\rho(g^{-1})\boldsymbol{\alpha}
 \qquad(g\in G).
\]
\end{definition}

\begin{proposition}[Forms with values in an associated bundle]
\label{prop:formas-horizontales-equivariantes-asociadas}
Let $M$ be a smooth manifold with or without boundary, and let
$\mathbf{E}=\mathbf{P}\times_\rho V\to M$ be an associated bundle. There is a natural isomorphism
\[
 \Omega^k(M;\mathbf{E})
 \cong
 \{\boldsymbol{\alpha}\in\Omega^k(\mathbf{P};V)\mid
 \boldsymbol{\alpha}\text{ is horizontal and }\rho\text{-equivariant}\}.
\]
If $\boldsymbol{\alpha}$ belongs to the second space, the associated form
$\underline{\boldsymbol{\alpha}}$ is defined by
\begin{equation}
\label{eq:descenso-forma-horizontal-equiv}
 \underline{\boldsymbol{\alpha}}_x(X_1,\dots,X_k)
 :=
 [p,\boldsymbol{\alpha}_p(\widetilde X_1,\dots,\widetilde X_k)]_\rho,
\end{equation}
where $p\in \mathbf{P}_x$ and $d\pi_p\widetilde X_i=X_i$.
\end{proposition}

\begin{proof}
Horizontality shows that
\eqref{eq:descenso-forma-horizontal-equiv} is independent of the lifts
$\widetilde X_i$. If $p$ is replaced by $p\cdot g$, we may choose the lifts $(dR_g)_p\widetilde X_i$; equivariance then gives
\[
 \boldsymbol{\alpha}_{p\cdot g}
 \bigl((dR_g)\widetilde X_1,\dots,(dR_g)\widetilde X_k\bigr)
 =
 \rho(g^{-1})\boldsymbol{\alpha}_p(\widetilde X_1,\dots,\widetilde X_k),
\]
and the equivalence relation of the associated bundle shows that the class in
\eqref{eq:descenso-forma-horizontal-equiv} is unchanged.

Conversely, given $\boldsymbol{\beta}\in\Omega^k(M;\mathbf{E})$ and $p\in \mathbf{P}_x$, there is a unique
$v\in V$ such that
\[
 \boldsymbol{\beta}_x(d\pi\,\widetilde X_1,\dots,d\pi\,\widetilde X_k)=[p,v]_\rho.
\]
Define $\boldsymbol{\alpha}_p(\widetilde X_1,\dots,\widetilde X_k):=v$. This form is horizontal, and uniqueness of $v$ proves equivariance. In a principal trivialization, both constructions reduce to the usual identifications with $V$-valued forms, proving smoothness and that the two operations are inverses.
\end{proof}

For $\xi\in\mathfrak g$, denote by
\[
 \boldsymbol{\xi}_{\mathbf{P}}(p):=
 \left.\frac{d}{dt}\right|_0p\cdot\exp_G(t\xi)
\]
the positive fundamental vector field of the right action. By
\eqref{eq:signo-campos-fundamentales-derecha},
\[
 [\boldsymbol{\xi}_{\mathbf{P}},\boldsymbol{\eta}_{\mathbf{P}}]
 =\boldsymbol{[\xi,\eta]}_{\mathbf{P}}.
\]
For each $p$, the map $\xi\mapsto\boldsymbol{\xi}_{\mathbf{P}}(p)$ is an isomorphism $\mathfrak g\to\mathbf{V}_p\mathbf{P}$.

\begin{definition}[Principal connection]
\label{def:forma-conexion-principal}
\index{principal connection}
Let $M$ be a smooth manifold with or without boundary, and let
$\mathbf{P}\to M$ be a right principal bundle. A \textit{principal connection form} on $\mathbf{P}$ is a form
$\boldsymbol{\omega}\in\Omega^1(\mathbf{P};\mathfrak g)$ such that
\begin{align}
 \boldsymbol{\omega}(\boldsymbol{\xi}_{\mathbf{P}})&=\xi,
 \label{eq:conexion-principal-reproduce}\\
 R_g^*\boldsymbol{\omega}&=\operatorname{Ad}_{g^{-1}}\boldsymbol{\omega}.
 \label{eq:conexion-principal-equiv}
\end{align}
The corresponding horizontal distribution is
\[
 \mathbf{H}_p\mathbf{P}:=\ker\boldsymbol{\omega}_p.
\]
Then
\[
 T_p\mathbf{P}=\mathbf{H}_p\mathbf{P}\oplus \mathbf{V}_p\mathbf{P},
 \qquad
 (dR_g)_p\mathbf{H}_p\mathbf{P}=\mathbf{H}_{p\cdot g}\mathbf{P}.
\]
\end{definition}

\begin{proposition}[Horizontal distributions and connection forms]
\label{prop:equivalencia-horizontal-forma-conexion-principal}
Let $M$ be a smooth manifold with or without boundary, and let
$\mathbf{P}\to M$ be a right principal bundle. The assignment
$\mathbf{H}_p\mathbf{P}=\ker\boldsymbol{\omega}_p$ establishes a bijection between principal connection forms and smooth distributions
$\mathbf{H}\subset T\mathbf{P}$ such that
\[
 T_p\mathbf{P}=\mathbf{H}_p\oplus \mathbf{V}_p\mathbf{P},
 \qquad
 (dR_g)_p\mathbf{H}_p=\mathbf{H}_{pg}.
\]
The form corresponding to $\mathbf{H}$ is obtained by writing uniquely
$Z=Z^{\mathbf{H}}+\boldsymbol{\xi}_{\mathbf{P}}(p)$ and setting
$\boldsymbol{\omega}_p(Z)=\xi$.
\end{proposition}

\begin{proof}
For a connection form, the first property implies that
$\boldsymbol{\omega}_p|_{\mathbf{V}_p\mathbf{P}}\colon
\mathbf{V}_p\mathbf{P}\to\mathfrak g$ is the inverse of
$\xi\mapsto\boldsymbol{\xi}_{\mathbf{P}}(p)$, giving the direct sum. The second property implies invariance of the kernels.

Conversely, the smooth isomorphism
$\mathfrak g\to \mathbf{V}_p\mathbf{P}$,
$\xi\mapsto\boldsymbol{\xi}_{\mathbf{P}}(p)$, and the smooth projection onto
$\mathbf{V}_p\mathbf{P}$ along $\mathbf{H}_p\mathbf{P}$ define the indicated form. We have
$\boldsymbol{\omega}(\boldsymbol{\xi}_{\mathbf{P}})=\xi$. Invariance of
$\mathbf{H}$ and the identity
$(dR_g)_p\boldsymbol{\xi}_{\mathbf{P}}(p)
=(\operatorname{Ad}_{g^{-1}}\xi)_{\mathbf{P}}(pg)$ imply
$R_g^*\boldsymbol{\omega}=\operatorname{Ad}_{g^{-1}}\boldsymbol{\omega}$. The two constructions are inverses by uniqueness of the horizontal--vertical decomposition.
\end{proof}

\begin{definition}[Wedge bracket]
\label{def:corchete-cuna-formas-lie}
\index{wedge bracket}
Let $N$ be a smooth manifold with or without boundary. For
$\boldsymbol{\alpha}\in\Omega^k(N;\mathfrak g)$ and
$\boldsymbol{\beta}\in\Omega^\ell(N;\mathfrak g)$, define
\[
 [\boldsymbol{\alpha}\wedge\boldsymbol{\beta}](\mathbf{X}_1,\dots,\mathbf{X}_{k+\ell})
 :=
 \displaystyle\sum_{\sigma\in\operatorname{Sh}(k,\ell)}
 \operatorname{sgn}(\sigma)
 \left[
 \boldsymbol{\alpha}(\mathbf{X}_{\sigma(1)},\dots,\mathbf{X}_{\sigma(k)}),
 \boldsymbol{\beta}(\mathbf{X}_{\sigma(k+1)},\dots,\mathbf{X}_{\sigma(k+\ell)})
 \right],
\]
where $\operatorname{Sh}(k,\ell)$ is the set of
$(k,\ell)$-shuffle permutations. In particular,
\[
 [\boldsymbol{\alpha}\wedge\boldsymbol{\beta}](\mathbf{X},\mathbf{Y})
 =
 [\boldsymbol{\alpha}(\mathbf{X}),\boldsymbol{\beta}(\mathbf{Y})]
 -[\boldsymbol{\alpha}(\mathbf{Y}),\boldsymbol{\beta}(\mathbf{X})]
\]
for $\boldsymbol{\alpha},\boldsymbol{\beta}\in\Omega^1(N;\mathfrak g)$.
\end{definition}

\begin{lemma}[Graded algebra of the wedge bracket]
\label{lem:algebra-graduada-corchete-cuna}
Let $N$ be a smooth manifold with or without boundary, and let
$\boldsymbol\alpha\in\Omega^k(N;\mathfrak g)$,
$\boldsymbol\beta\in\Omega^\ell(N;\mathfrak g)$, and
$\boldsymbol\gamma\in\Omega^m(N;\mathfrak g)$, with
$k,\ell,m\in\mathbb N_0$. Then
\begin{align}
[\boldsymbol\beta\wedge\boldsymbol\alpha]
&=-(-1)^{k\ell}[\boldsymbol\alpha\wedge\boldsymbol\beta],
\label{eq:antisimetria-graduada-corchete-cuna}\\
d[\boldsymbol\alpha\wedge\boldsymbol\beta]
&=[d\boldsymbol\alpha\wedge\boldsymbol\beta]
+(-1)^k[\boldsymbol\alpha\wedge d\boldsymbol\beta].
\label{eq:leibniz-corchete-cuna}
\end{align}
The graded Jacobi identity is
\[
\begin{aligned}
&(-1)^{km}[\boldsymbol\alpha\wedge[\boldsymbol\beta\wedge\boldsymbol\gamma]]
+(-1)^{\ell k}[\boldsymbol\beta\wedge[\boldsymbol\gamma\wedge\boldsymbol\alpha]]\\
&\qquad+(-1)^{m\ell}[\boldsymbol\gamma\wedge[\boldsymbol\alpha\wedge\boldsymbol\beta]]=0.
\end{aligned}
\]
In particular, for $\mathbf A\in\Omega^1(N;\mathfrak g)$,
\begin{equation}
\label{eq:jacobi-A-A-A}
[\mathbf A\wedge[\mathbf A\wedge\mathbf A]]=0.
\end{equation}
\end{lemma}
\begin{proof}
Fix a basis $(e_1,\ldots,e_d)$ of $\mathfrak g$. Write
$\boldsymbol\alpha=\displaystyle\sum_{a=1}^d\alpha^a\otimes e_a$ and
$\boldsymbol\beta=\displaystyle\sum_{b=1}^d\beta^b\otimes e_b$,
where $\alpha^a,\beta^b$ are scalar forms. The definition using
shuffles is equivalent to
\[
[\boldsymbol\alpha\wedge\boldsymbol\beta]
=\displaystyle\sum_{a,b=1}^d(\alpha^a\wedge\beta^b)\otimes[e_a,e_b].
\]
Interchanging the scalar forms produces $(-1)^{k\ell}$, while
interchanging the arguments of the Lie bracket produces $-1$.
This proves the first identity. Apply $d$ to each product
$\alpha^a\wedge\beta^b$ and use its Leibniz rule to prove the second.

Write also
$\boldsymbol\gamma=\displaystyle\sum_{c=1}^d\gamma^c\otimes e_c$.
In all three terms of graded Jacobi, order the scalar forms as
$\alpha^a\wedge\beta^b\wedge\gamma^c$. The additional sign in the
second term is $(-1)^{k(\ell+m)}$; in the third it is
$(-1)^{m(k+\ell)}$. Multiplying by the signs in the statement makes
all three coefficients $(-1)^{km}$. The sum becomes
\[
\begin{aligned}
(-1)^{km}\displaystyle\sum_{a,b,c=1}^d
&(\alpha^a\wedge\beta^b\wedge\gamma^c)\otimes\\
&\bigl([e_a,[e_b,e_c]]+[e_b,[e_c,e_a]]+[e_c,[e_a,e_b]]\bigr)=0,
\end{aligned}
\]
by Jacobi in $\mathfrak g$. If all three forms are $\mathbf A$ of
degree one, all three summands are
$-[\mathbf A\wedge[\mathbf A\wedge\mathbf A]]$. Dividing by $-3$
proves \eqref{eq:jacobi-A-A-A}.
\end{proof}

\begin{definition}[Principal curvature]
\label{def:curvatura-conexion-principal}
Let $M$ be a smooth manifold with or without boundary, and let $\boldsymbol{\omega}$ be a principal connection on $\mathbf{P}\to M$. The curvature of
$\boldsymbol{\omega}$ is
\begin{equation}
\label{eq:curvatura-principal}
 \boldsymbol{\Omega}_{\boldsymbol{\omega}}
 :=d\boldsymbol{\omega}+\frac12[\boldsymbol{\omega}\wedge\boldsymbol{\omega}]
 \in\Omega^2(\mathbf{P};\mathfrak g).
\end{equation}
\end{definition}

We will use an elementary identity for forms and flows.

\begin{lemma}[Cartan's formula]
\label{lem:formula-cartan-flujo}
Let $N$ be a smooth manifold with or without boundary, let $\mathbf{X}$ be a vector field whose local flow remains in $N$, and let $\boldsymbol{\alpha}$ be a differential form. Then
\[
 \mathcal L_{\mathbf{X}}\boldsymbol{\alpha}
 :=
 \left.\frac{d}{dt}\right|_0\Phi_t^*\boldsymbol{\alpha}
 =
 d(\iota_{\mathbf{X}}\boldsymbol{\alpha})
 +\iota_{\mathbf{X}}(d\boldsymbol{\alpha}).
\]
\end{lemma}

\begin{proof}
Both sides are graded derivations of degree zero for the exterior product and commute with restrictions. For a function $f$, both sides equal $\mathbf{X}f$. In coordinates, every form is a sum of products of functions and differentials $\mathbf{d}x^i$; on the latter,
\[
 \mathcal L_{\mathbf{X}}(\mathbf{d}x^i)=d(\mathbf{X}x^i)=d(\iota_{\mathbf{X}}\mathbf{d}x^i).
\]
The Leibniz rule establishes the identity for every form.
\end{proof}

\begin{proposition}
\label{prop:curvatura-horizontal-equiv}
Let $M$ be a smooth manifold with or without boundary, and let $\boldsymbol{\omega}$ be a principal connection on $\mathbf{P}\to M$. The curvature
$\boldsymbol{\Omega}_{\boldsymbol{\omega}}$ is horizontal and
$\operatorname{Ad}$-equivariant:
\[
 R_g^*\boldsymbol{\Omega}_{\boldsymbol{\omega}}
 =
 \operatorname{Ad}_{g^{-1}}\boldsymbol{\Omega}_{\boldsymbol{\omega}}.
\]
Thus, by
Proposition~\ref{prop:formas-horizontales-equivariantes-asociadas},
it determines a global form
\[
 \mathbf{F}_{\boldsymbol{\omega}}\in\Omega^2(M;\operatorname{ad}(\mathbf{P})).
\]
\end{proposition}

\begin{proof}
Equivariance follows from
$R_g^*d\boldsymbol{\omega}=d(R_g^*\boldsymbol{\omega})$ and from the fact that
$\operatorname{Ad}_{g^{-1}}$ preserves the bracket.

The flow of $\boldsymbol{\xi}_{\mathbf{P}}$ is $R_{\exp(t\xi)}$.
Differentiating the equivariance of
$\boldsymbol{\omega}$ gives
\[
 \mathcal L_{\boldsymbol{\xi}_{\mathbf{P}}}\boldsymbol{\omega}
 =
 -\operatorname{ad}_\xi\boldsymbol{\omega}
 =
 -[\xi,\boldsymbol{\omega}].
\]
Since $\iota_{\boldsymbol{\xi}_{\mathbf{P}}}\boldsymbol{\omega}=\xi$ is constant, Cartan's formula gives
\[
 \iota_{\boldsymbol{\xi}_{\mathbf{P}}}d\boldsymbol{\omega}
 =-[\xi,\boldsymbol{\omega}].
\]
On the other hand, the definition of the wedge bracket implies
\[
 \frac12\iota_{\boldsymbol{\xi}_{\mathbf{P}}}
 [\boldsymbol{\omega}\wedge\boldsymbol{\omega}]=[\xi,\boldsymbol{\omega}].
\]
The two terms cancel, so
$\iota_{\boldsymbol{\xi}_{\mathbf{P}}}
\boldsymbol{\Omega}_{\boldsymbol{\omega}}=0$. Every vertical vector has the form $\boldsymbol{\xi}_{\mathbf{P}}(p)$, and hence the curvature is horizontal.
\end{proof}

\begin{proposition}[Local potentials and curvatures]
\label{prop:potenciales-curvaturas-locales-principal}
Let $M$ be a smooth manifold with or without boundary, let $\boldsymbol{\omega}$ be a principal connection on $\mathbf{P}\to M$, and let
$\mathbf{s}_\alpha\colon U_\alpha\to \mathbf{P}$ be a local section.
Define
\[
 \mathbf{A}_\alpha:=\mathbf{s}_\alpha^*\boldsymbol{\omega}
 \in\Omega^1(U_\alpha;\mathfrak g),
 \qquad
 \mathbf{F}_\alpha:=\mathbf{s}_\alpha^*
 \boldsymbol{\Omega}_{\boldsymbol{\omega}}.
\]
Then
\begin{equation}
\label{eq:curvatura-local-conexion-principal}
 \mathbf{F}_\alpha=d\mathbf{A}_\alpha
 +\frac12[\mathbf{A}_\alpha\wedge \mathbf{A}_\alpha].
\end{equation}
If $\mathbf{s}_\beta=\mathbf{s}_\alpha g_{\alpha\beta}$, then
\begin{align}
 \mathbf{A}_\beta
 &=
 \operatorname{Ad}_{g_{\alpha\beta}^{-1}}\mathbf{A}_\alpha
 +g_{\alpha\beta}^{-1}dg_{\alpha\beta},
 \label{eq:transformacion-potencial-conexion-principal}\\
 \mathbf{F}_\beta
 &=
 \operatorname{Ad}_{g_{\alpha\beta}^{-1}}\mathbf{F}_\alpha.
 \label{eq:transformacion-curvatura-principal}
\end{align}
Here $g^{-1}dg$ denotes the pullback of the left Maurer--Cartan form of $G$.
\end{proposition}

\begin{proof}
The first identity is the pullback of
\eqref{eq:curvatura-principal}. Let $X\in T_xM$ and write
$g=g_{\alpha\beta}(x)$. Differentiating
$\mathbf{s}_\beta=\mathbf{s}_\alpha g_{\alpha\beta}$ gives
\[
 d(\mathbf{s}_\beta)_xX
 =
 (dR_g)_{\mathbf{s}_\alpha(x)}d(\mathbf{s}_\alpha)_xX
 +
 \bigl(g^{-1}dg(X)\bigr)_{\mathbf{P}}
 \big|_{\mathbf{s}_\beta(x)}.
\]
Applying $\boldsymbol{\omega}$ and using
\eqref{eq:conexion-principal-reproduce}--\eqref{eq:conexion-principal-equiv}
yields
\eqref{eq:transformacion-potencial-conexion-principal}.

For curvature, the second summand of the preceding expression is vertical, while $\boldsymbol{\Omega}_{\boldsymbol{\omega}}$ is horizontal. The horizontal summands transform by $dR_g$, and equivariance of
$\boldsymbol{\Omega}_{\boldsymbol{\omega}}$ gives
\eqref{eq:transformacion-curvatura-principal}.
\end{proof}

\begin{definition}[Exterior covariant derivative]
\label{def:derivada-exterior-covariante-principal}
Let $M$ be a smooth manifold with or without boundary, let $\boldsymbol{\omega}$ be a principal connection on $\mathbf{P}\to M$, let
$\rho\colon G\to GL(V)$ be a representation, and let
$\rho_*\colon\mathfrak g\to\operatorname{End}(V)$ be its differential. For a horizontal,
$\rho$-equivariant form
$\boldsymbol{\alpha}\in\Omega^k(\mathbf{P};V)$, define
\begin{equation}
\label{eq:Domega-forma-equiv}
 D_{\boldsymbol{\omega}}\boldsymbol{\alpha}
 :=
 d\boldsymbol{\alpha}+\rho_*(\boldsymbol{\omega})\wedge\boldsymbol{\alpha}.
\end{equation}
The second operation is given by
\[
 \bigl(\rho_*(\boldsymbol{\omega})\wedge\boldsymbol{\alpha}\bigr)
 (\mathbf{X}_0,\dots,\mathbf{X}_k)
 =
 \displaystyle\sum_{i=0}^k(-1)^i
 \rho_*(\boldsymbol{\omega}(\mathbf{X}_i))
 \boldsymbol{\alpha}(\mathbf{X}_0,\dots,
 \widehat{\mathbf{X}}_i,\dots,\mathbf{X}_k).
\]
\end{definition}

\begin{proposition}
\label{prop:Domega-desciende-dA}
Let $M$ be a smooth manifold with or without boundary, and let
$\boldsymbol{\omega}$, $\rho$, and $\boldsymbol{\alpha}$ be as in the preceding definition. The form $D_{\boldsymbol{\omega}}\boldsymbol{\alpha}$ is again horizontal and
$\rho$-equivariant. Consequently, it descends to an operator
\[
 d_{\boldsymbol{\omega}}\colon
 \Omega^k(M;\mathbf{P}\times_\rho V)
 \longrightarrow
 \Omega^{k+1}(M;\mathbf{P}\times_\rho V).
\]
In a local section $\mathbf{s}_\alpha$, if
$\mathbf{a}_\alpha:=\mathbf{s}_\alpha^*\boldsymbol{\alpha}$, then
\begin{equation}
\label{eq:dA-local-representacion}
 (d_{\boldsymbol{\omega}}\mathbf{a})_\alpha
 =
 d_{\mathbf{A}_\alpha}\mathbf{a}_\alpha
 :=d\mathbf{a}_\alpha
 +\rho_*(\mathbf{A}_\alpha)\wedge \mathbf{a}_\alpha.
\end{equation}
If $\mathbf{a}_\beta=\rho(g_{\alpha\beta}^{-1})\mathbf{a}_\alpha$, then
\begin{equation}
\label{eq:covarianza-dA-local}
 d_{\mathbf{A}_\beta}\mathbf{a}_\beta
 =
 \rho(g_{\alpha\beta}^{-1})d_{\mathbf{A}_\alpha}\mathbf{a}_\alpha.
\end{equation}
For the adjoint representation,
\[
 d_{\mathbf{A}}\mathbf{a}
 =d\mathbf{a}+[\mathbf{A}\wedge \mathbf{a}].
\]
\end{proposition}

\begin{proof}
Equivariance of $D_{\boldsymbol{\omega}}\boldsymbol{\alpha}$ follows by applying
$R_g^*$ to \eqref{eq:Domega-forma-equiv} and using
\[
 \rho_*(\operatorname{Ad}_{g^{-1}}\xi)
 =
 \rho(g^{-1})\rho_*(\xi)\rho(g).
\]
To prove horizontality, let $\boldsymbol{\xi}_{\mathbf{P}}$ be a fundamental vertical vector field. Infinitesimal equivariance of
$\boldsymbol{\alpha}$ gives
\[
 \mathcal L_{\boldsymbol{\xi}_{\mathbf{P}}}\boldsymbol{\alpha}
 =-\rho_*(\xi)\boldsymbol{\alpha}.
\]
Since $\iota_{\boldsymbol{\xi}_{\mathbf{P}}}\boldsymbol{\alpha}=0$, Cartan's formula implies
\[
 \iota_{\boldsymbol{\xi}_{\mathbf{P}}}d\boldsymbol{\alpha}
 =-\rho_*(\xi)\boldsymbol{\alpha}.
\]
On the other hand,
\[
 \iota_{\boldsymbol{\xi}_{\mathbf{P}}}
 \bigl(\rho_*(\boldsymbol{\omega})\wedge\boldsymbol{\alpha}\bigr)
 =
 \rho_*(\xi)\boldsymbol{\alpha},
\]
and the terms cancel. Thus
$D_{\boldsymbol{\omega}}\boldsymbol{\alpha}$ descends by
Proposition~\ref{prop:formas-horizontales-equivariantes-asociadas}.

The local formula follows by pulling back by $\mathbf{s}_\alpha$.
The covariance identity expresses that two local pullbacks represent the same global form. It can also be checked directly by substituting
\eqref{eq:transformacion-potencial-conexion-principal} and
$\mathbf{a}_\beta=\rho(g_{\alpha\beta}^{-1})\mathbf{a}_\alpha$. For
$\rho=\operatorname{Ad}$, we have
$\rho_*(\xi)\eta=[\xi,\eta]$.
\end{proof}

\begin{corollary}[Induced connection on an associated bundle]
\label{cor:conexion-inducida-haz-asociado-principal}
Let $M$ be a smooth manifold with or without boundary, let $\boldsymbol{\omega}$ be a principal connection on $\mathbf{P}\to M$, and let
$\rho\colon G\to GL(V)$ be a representation. For
$\mathbf{E}=\mathbf{P}\times_\rho V$, the degree-zero operator
\[
 \nabla^{\mathbf{E},\boldsymbol{\omega}}:=d_{\boldsymbol{\omega}}\colon
 \Gamma(\mathbf{E})=\Omega^0(M;\mathbf{E})\longrightarrow\Omega^1(M;\mathbf{E})
\]
is a complex or real connection according to $V$. If $u_\alpha$ are the local components of a section $\mathbf{u}$, then
\begin{equation}
\label{eq:conexion-inducida-asociada-local}
 (\nabla^{\mathbf{E},\boldsymbol{\omega}}\mathbf{u})_\alpha
 =du_\alpha+\rho_*(\mathbf{A}_\alpha)u_\alpha.
\end{equation}
Equivalently, if $\widetilde u\colon \mathbf{P}\to V$ is the associated equivariant map and $\widetilde{\mathbf{X}}$ is the horizontal lift of a vector field
$\mathbf{X}$, then
\[
 (\nabla_{\mathbf{X}}^{\mathbf{E},\boldsymbol{\omega}}\mathbf{u})(x)
 =[p,d\widetilde u_p(\widetilde{\mathbf{X}}_p)]_\rho.
\]
\end{corollary}

\begin{proof}
The local formula \eqref{eq:dA-local-representacion} with $k=0$ is
\eqref{eq:conexion-inducida-asociada-local}. It gives linearity and
$\nabla(f\mathbf{u})=df\otimes \mathbf{u}+f\nabla \mathbf{u}$. The covariance identity
\eqref{eq:covarianza-dA-local} shows that the expressions glue. In the equivariant model, $\boldsymbol{\omega}(\widetilde{\mathbf{X}})=0$, so
$D_{\boldsymbol{\omega}}\widetilde u(\widetilde{\mathbf{X}})
=d\widetilde u(\widetilde{\mathbf{X}})$; the descent of forms proves the final formula.
\end{proof}

\begin{theorem}[Square of the covariant derivative and the Bianchi identity]
\label{teo:dA-cuadrado-bianchi-principal}
Let $M$ be a smooth manifold with or without boundary, and let $\mathbf{A}$ be a local potential for a principal connection on $\mathbf{P}\to M$. For every representation $\rho$ and every local form $\mathbf{a}$ with values in $V$,
\begin{equation}
\label{eq:dA-cuadrado-curvatura}
 d_{\mathbf{A}}^2\mathbf{a}
 =\rho_*(\mathbf{F}_{\mathbf{A}})\wedge \mathbf{a}.
\end{equation}
In particular, for forms with values in $\operatorname{ad}(\mathbf{P})$,
\begin{equation}
\label{eq:dA-cuadrado-adjunto}
 d_{\mathbf{A}}^2\mathbf{a}=[\mathbf{F}_{\mathbf{A}}\wedge \mathbf{a}].
\end{equation}
The curvature satisfies the Bianchi identity
\begin{equation}
\label{eq:bianchi-local-principal}
 d_{\mathbf{A}}\mathbf{F}_{\mathbf{A}}
 =d\mathbf{F}_{\mathbf{A}}
 +[\mathbf{A}\wedge \mathbf{F}_{\mathbf{A}}]=0.
\end{equation}
These identities are compatible with changes of local section and are therefore global.
\end{theorem}

\begin{proof}
Write $\mathbf{B}=\rho_*(\mathbf{A})$. Using the exterior Leibniz rule,
\[
 (d+\mathbf{B}\wedge)^2\mathbf{a}
 =
 d\mathbf{B}\wedge \mathbf{a}
 +\mathbf{B}\wedge \mathbf{B}\wedge \mathbf{a}.
\]
Since $\rho_*$ is a Lie algebra homomorphism,
\[
 \mathbf{B}\wedge \mathbf{B}
 =
 \frac12\rho_*([\mathbf{A}\wedge \mathbf{A}]).
\]
This proves \eqref{eq:dA-cuadrado-curvatura}, and the adjoint representation gives \eqref{eq:dA-cuadrado-adjunto}.

For the Bianchi identity, by
$\mathbf{F}_{\mathbf{A}}=d\mathbf{A}
+\frac12[\mathbf{A}\wedge \mathbf{A}]$ and
Lemma~\ref{lem:algebra-graduada-corchete-cuna},
\[
 d\mathbf{F}_{\mathbf{A}}
 =
 \frac12\bigl([d\mathbf{A}\wedge \mathbf{A}]
 -[\mathbf{A}\wedge d\mathbf{A}]\bigr)
 =
 -[\mathbf{A}\wedge d\mathbf{A}].
\]
Moreover,
\[
 [\mathbf{A}\wedge \mathbf{F}_{\mathbf{A}}]
 =
 [\mathbf{A}\wedge d\mathbf{A}]
 +\frac12[\mathbf{A}\wedge[\mathbf{A}\wedge \mathbf{A}]]
 =
 [\mathbf{A}\wedge d\mathbf{A}].
\]
The sum is zero. The covariance identities
\eqref{eq:covarianza-dA-local} and
\eqref{eq:transformacion-curvatura-principal} show that all these identities glue.
\end{proof}

\begin{theorem}[Existence of principal connections]
\label{teo:existencia-conexiones-principales}
Let $M$ be a smooth manifold with or without boundary. Every smooth principal bundle
$\mathbf{P}\to M$ admits a principal connection.
\end{theorem}

\begin{proof}
Let $\{U_\alpha\}_{\alpha\in I}$ be a trivializing cover. On
$\pi^{-1}(U_\alpha)\cong U_\alpha\times G$, define
\[
 \boldsymbol{\omega}_\alpha(X,\zeta):=\theta^L_g(\zeta),
 \qquad
 (X,\zeta)\in T_xU_\alpha\oplus T_gG,
\]
where $\theta^L_g=(dL_{g^{-1}})_g$ is the left Maurer--Cartan form.
A direct verification gives
\[
 \boldsymbol{\omega}_\alpha(\boldsymbol{\xi}_{\mathbf{P}})=\xi,\qquad
 R_h^*\boldsymbol{\omega}_\alpha=\operatorname{Ad}_{h^{-1}}\boldsymbol{\omega}_\alpha;
\]
therefore $\boldsymbol{\omega}_\alpha$ is a local principal connection.

Let $\{\chi_\alpha\}_{\alpha\in I}$ be a smooth partition of unity subordinate to
$\{U_\alpha\}$, whose existence is given by
Theorem~\ref{particionesdelaunidad}. Each form
$(\pi^*\chi_\alpha)\boldsymbol{\omega}_\alpha$, extended by zero outside
$\pi^{-1}(U_\alpha)$, is smooth, and the family is locally finite. Define
\[
 \boldsymbol{\omega}:=\displaystyle\sum_{\alpha\in I}(\pi^*\chi_\alpha)\boldsymbol{\omega}_\alpha.
\]
Since $\displaystyle \sum_{\alpha\in I}\chi_\alpha=1$,
\[
 \boldsymbol{\omega}(\boldsymbol{\xi}_{\mathbf{P}})
 =\displaystyle\sum_{\alpha\in I}(\pi^*\chi_\alpha)\xi=\xi.
\]
The functions $\pi^*\chi_\alpha$ are invariant under the right action and each $\boldsymbol{\omega}_\alpha$ is equivariant, so
$R_g^*\boldsymbol{\omega}=\operatorname{Ad}_{g^{-1}}\boldsymbol{\omega}$. Therefore $\boldsymbol{\omega}$ is a global principal connection.
\end{proof}

\begin{theorem}[Affine structure of the space of connections]
\label{teo:espacio-afin-conexiones-principales}
Let $M$ be a smooth manifold with or without boundary, and let
$\mathbf{P}\to M$ be a principal bundle. Let $\operatorname{Conn}(\mathbf{P})$ be the set of principal connections on
$\mathbf{P}$. If $\boldsymbol{\omega}_0,\boldsymbol{\omega}_1\in\operatorname{Conn}(\mathbf{P})$, then
\[
 \mathbf{a}:=\boldsymbol{\omega}_1-\boldsymbol{\omega}_0
\]
is horizontal and $\operatorname{Ad}$-equivariant, and therefore corresponds to a unique element of
\[
 \Omega^1(M;\operatorname{ad}(\mathbf{P})).
\]
Conversely, if $\boldsymbol{\omega}_0$ is a connection and
$\mathbf{a}\in\Omega^1(M;\operatorname{ad}(\mathbf{P}))$, the corresponding horizontal equivariant form
$\widetilde{\mathbf{a}}$ satisfies the condition that
\[
 \boldsymbol{\omega}_0+\widetilde{\mathbf{a}}
\]
is a principal connection. Consequently,
$\operatorname{Conn}(\mathbf{P})$ is an affine space modeled on
$\Omega^1(M;\operatorname{ad}(\mathbf{P}))$.
\end{theorem}

\begin{proof}
For every $\xi\in\mathfrak g$,
\[
 \mathbf{a}(\boldsymbol{\xi}_{\mathbf{P}})
 =\boldsymbol{\omega}_1(\boldsymbol{\xi}_{\mathbf{P}})
 -\boldsymbol{\omega}_0(\boldsymbol{\xi}_{\mathbf{P}})=\xi-\xi=0,
\]
so $\mathbf{a}$ is horizontal. Equivariance follows by subtracting the equivariance identities for $\boldsymbol{\omega}_1$ and $\boldsymbol{\omega}_0$.
Proposition~\ref{prop:formas-horizontales-equivariantes-asociadas} then identifies $\mathbf{a}$ with a form in
$\Omega^1(M;\operatorname{ad}(\mathbf{P}))$.

Conversely, a horizontal equivariant form vanishes on fundamental vector fields and satisfies the same transformation law as a connection form. Therefore,
\[
 (\boldsymbol{\omega}_0+\widetilde{\mathbf{a}})
 (\boldsymbol{\xi}_{\mathbf{P}})=\xi,\qquad
 R_g^*(\boldsymbol{\omega}_0+\widetilde{\mathbf{a}})
 =
 \operatorname{Ad}_{g^{-1}}(\boldsymbol{\omega}_0+\widetilde{\mathbf{a}}).
\]
\end{proof}

\begin{corollary}
\label{cor:metrica-adP-grupo-compacto}
Let $M$ be a smooth manifold with or without boundary, and let
$\mathbf{P}\to M$ be a principal bundle with compact group $G$. Any
$\operatorname{Ad}$-invariant inner product given by
Corollary~\ref{cor:producto-interno-Ad-invariante-compacto} induces a well-defined bundle metric on $\operatorname{ad}(\mathbf{P})$ by
\[
 \langle[p,\xi],[p,\eta]\rangle_{\operatorname{ad}(\mathbf{P})}
 :=
 \langle\xi,\eta\rangle_{\mathfrak g}.
\]
\end{corollary}

\begin{proof}
Replacing $p$ by $p\cdot g$ replaces the representatives by
$\operatorname{Ad}_{g^{-1}}\xi$ and
$\operatorname{Ad}_{g^{-1}}\eta$. Invariance of the inner product shows that the value does not change.
\end{proof}

\subsection{Parallel transport and holonomy}

\begin{theorem}[Horizontal lift]
\label{teo:elevacion-horizontal-principal-unica}
\index{horizontal lift}
Let $M$ be a smooth manifold with or without boundary, let $\boldsymbol{\omega}$ be a principal connection on $\mathbf{P}\to M$, let
$\gamma\colon[a,b]\to M$ be a piecewise smooth curve, and let
$p_0\in \mathbf{P}_{\gamma(a)}$. There exists a unique piecewise smooth curve
$\widetilde\gamma_{p_0}\colon[a,b]\to \mathbf{P}$ such that
\[
 \widetilde\gamma_{p_0}(a)=p_0,\qquad
 \pi\circ\widetilde\gamma_{p_0}=\gamma,\qquad
 \dot{\widetilde\gamma}_{p_0}(t)
 \in \mathbf{H}_{\widetilde\gamma_{p_0}(t)}\mathbf{P}
\]
on each interval of smoothness. Moreover,
\[
 \widetilde\gamma_{p_0g}(t)
 =\widetilde\gamma_{p_0}(t)g.
\]
\end{theorem}
\begin{proof}
In a trivialization, write
$\widetilde\gamma(t)=\mathbf s(\gamma(t))g(t)$ and
$\mathbf A=\mathbf s^*\boldsymbol\omega$. Evaluating the connection
form on the derivative of this product shows that horizontality is
equivalent to
\[
\theta^L_{g(t)}\dot g(t)
=-\operatorname{Ad}_{g(t)^{-1}}
\bigl(\mathbf A_{\gamma(t)}(\dot\gamma(t))\bigr).
\]
If $a(t):=\mathbf A_{\gamma(t)}(\dot\gamma(t))\in\mathfrak g$,
this equation is equivalently
\[
\dot g(t)=-(dR_{g(t)})_e a(t).
\]
Indeed, $(dL_{g^{-1}})_g(dR_g)_e=\operatorname{Ad}_{g^{-1}}$.
This is a smooth time-dependent differential equation on $G$.
Local existence alone does not imply that its solution covers the
whole curve segment; we prove this extension using right invariance.

Fix a compact subinterval $[u,v]$ on which $\gamma$ is smooth and
stays in this trivialization. The coefficient $a(t)$ is smooth up
to the endpoints. In a chart of $G$ around $e$, local existence
provides a time $\delta>0$, uniform for every initial time
$\tau\in[u,v]$, for solutions $h_\tau$ with $h_\tau(\tau)=e$.
To justify uniformity, apply local existence with the initial time
as a parameter and cover the compact set $[u,v]\times\{e\}$ by
finitely many existence neighborhoods, reducing their times to the
positive minimum. Only values within $[u,v]$ are used.

For an initial value $g_0$ at time $\tau$, the curve
$h_\tau(t)g_0$ solves the same equation with that initial value,
because
\[
\frac{d}{dt}(h_\tau(t)g_0)
=-(dR_{g_0})_{h_\tau(t)}(dR_{h_\tau(t)})_e a(t)
=-(dR_{h_\tau(t)g_0})_e a(t).
\]
Thus the same $\delta$ works for every $g_0\in G$. Subdivide
$[u,v]$ into intervals of length less than $\delta$ and use each
final value as the next initial value. This produces the unique
solution on all of $[u,v]$. Local uniqueness guarantees agreement
on common domains.

Compactness of $[a,b]$ allows us to refine the smoothness partition
of $\gamma$ so that each segment lies in a trivializing domain.
Apply the preceding construction successively, starting each
segment at the endpoint of the previous lift. The lifts join to a
continuous piecewise smooth curve and are unique by uniqueness in
each trivialization. Finally, equivariance of the horizontal
distribution under $R_g$ shows that $\widetilde\gamma_{p_0}g$ is
horizontal and starts at $p_0g$. Uniqueness gives
$\widetilde\gamma_{p_0g}=\widetilde\gamma_{p_0}g$.
\end{proof}

\begin{definition}[Principal parallel transport]
\label{def:transporte-paralelo-principal}
\index{parallel transport!principal}
Let $M$ be a smooth manifold with or without boundary, let $\boldsymbol{\omega}$ be a principal connection on $\mathbf{P}\to M$, and let
$\gamma\colon[a,b]\to M$ be a piecewise smooth curve. The
\textbf{principal parallel transport} of $\gamma$ is
\[
 \mathcal P^{\boldsymbol{\omega}}_\gamma\colon
 \mathbf{P}_{\gamma(a)}\longrightarrow \mathbf{P}_{\gamma(b)},\qquad
 \mathcal P^{\boldsymbol{\omega}}_\gamma(p)
 :=\widetilde\gamma_p(b).
\]
\end{definition}

\begin{proposition}[Laws of principal transport]
\label{prop:leyes-transporte-paralelo-principal}
Let $M$ be a smooth manifold with or without boundary, and let $\boldsymbol{\omega}$ be a principal connection on $\mathbf{P}\to M$. Principal transport is
$G$-equivariant and satisfies
\[
 \mathcal P^{\boldsymbol{\omega}}_{\gamma^{-1}}
 =(\mathcal P^{\boldsymbol{\omega}}_\gamma)^{-1},\qquad
 \mathcal P^{\boldsymbol{\omega}}_{\gamma_2*\gamma_1}
 =\mathcal P^{\boldsymbol{\omega}}_{\gamma_2}\circ
  \mathcal P^{\boldsymbol{\omega}}_{\gamma_1},
\]
where $\gamma_2*\gamma_1$ traverses $\gamma_1$ first and then
$\gamma_2$. It is also invariant under orientation-preserving reparametrizations.
\end{proposition}
\begin{proof}
Equivariance is the final assertion of the preceding theorem. Reversing the parameter of a horizontal lift produces a horizontal lift of the reversed curve. For a concatenation, glue the lift of
$\gamma_1$ to the lift of $\gamma_2$ whose initial point is the endpoint of the first; uniqueness gives the formula. The same argument proves invariance under reparametrization.
\end{proof}

\begin{definition}[Principal holonomy group]
\label{def:grupo-holonomia-principal}
\index{holonomy group!principal}
Let $M$ be a smooth manifold with or without boundary, let $\boldsymbol{\omega}$ be a principal connection on $\mathbf{P}\to M$, and let $x=\pi(p)$. The
\textbf{holonomy group of $\boldsymbol{\omega}$ based at $p$} is
\[
 \operatorname{Hol}_p(\boldsymbol{\omega})
 :=
 \{g\in G\mid
   \mathcal P^{\boldsymbol{\omega}}_\gamma(p)=p g
   \text{ for some piecewise smooth loop }\gamma
   \text{ based at }x\}.
\]
\end{definition}

\begin{proposition}[Structure and change of base point for holonomy]
\label{prop:holonomia-principal-subgrupo-conjugacion}
Let $M$ be a smooth manifold with or without boundary, let $\boldsymbol{\omega}$ be a principal connection on $\mathbf{P}\to M$, and let $p\in\mathbf{P}$. The set $\operatorname{Hol}_p(\boldsymbol{\omega})$ is a subgroup of $G$.
If
$h\in G$, then
\[
 \operatorname{Hol}_{ph}(\boldsymbol{\omega})
 =h^{-1}\operatorname{Hol}_p(\boldsymbol{\omega})h.
\]
If $\lambda$ is a curve from $x$ to $y$ and
$q=\mathcal P^{\boldsymbol{\omega}}_\lambda(p)$, then
$\operatorname{Hol}_q(\boldsymbol{\omega})
=\operatorname{Hol}_p(\boldsymbol{\omega})$. For another choice $q h\in \mathbf{P}_y$, the groups are conjugate by $h$.
\end{proposition}
\begin{proof}
The constant loop represents $e$. If the loops $\alpha,\beta$ represent $g,h$, respectively, traversing $\beta$ and then
$\alpha$ takes us, by equivariance, from $p$ to $pgh$; reversing a loop representing $g$ produces $g^{-1}$. This proves the subgroup property.
Starting at $ph$, the lift of a loop that starts at $p$ and ends at
$pg$ instead ends at $pgh=(ph)(h^{-1}gh)$, giving the conjugation formula.
Finally, conjugating a loop based at $y$ by $\lambda$ and
$\lambda^{-1}$ preserves the element of $G$ when $q$ is the transport of $p$. The final assertion combines this fact with the formula within a single fiber.
\end{proof}

\begin{proposition}[Holonomy of an associated bundle]
\label{prop:holonomia-haz-vectorial-asociado}
\index{holonomy group!of a vector bundle}
Let $M$ be a smooth manifold with or without boundary, let $\boldsymbol{\omega}$ be a principal connection on $\mathbf{P}\to M$, and let
$\mathbf{E}=\mathbf{P}\times_\rho V$. The transport induced by
$\boldsymbol{\omega}$ is
\[
 \mathcal P^{\mathbf{E},\boldsymbol{\omega}}_\gamma[p,v]
 :=
 [\mathcal P^{\boldsymbol{\omega}}_\gamma(p),v].
\]
It is well defined, linear, and satisfies the inversion and concatenation laws. Therefore,
\[
 \operatorname{Hol}_x(\mathbf{E},\boldsymbol{\omega})
 :=\{\mathcal P^{\mathbf{E},\boldsymbol{\omega}}_\gamma\mid
       \gamma\text{ is a loop based at }x\}
 \subseteq GL(\mathbf{E}_x)
\]
is a subgroup. If
$\mathcal P^{\boldsymbol{\omega}}_\gamma(p)=pg$, the identification
$\mathbf{E}_x\cong V$ determined by $p$ represents
$\mathcal P^{\mathbf{E},\boldsymbol{\omega}}_\gamma$ by $\rho(g)$. When the base point changes, the vector bundle holonomy groups are conjugate through parallel transport along a curve joining the points.
\end{proposition}
\begin{proof}
If $[pg,v]=[p,\rho(g)v]$, principal equivariance gives
\[
 [\mathcal P^{\boldsymbol{\omega}}_\gamma(pg),v]
 =[\mathcal P^{\boldsymbol{\omega}}_\gamma(p)g,v]
 =[\mathcal P^{\boldsymbol{\omega}}_\gamma(p),\rho(g)v],
\]
so the definition is independent of the representative. The remaining assertions follow from
Proposition~\ref{prop:leyes-transporte-paralelo-principal}.
\end{proof}

\begin{definition}[Parallel section of the adjoint group bundle]
\label{def:seccion-paralela-AdP}
Let $M$ be a smooth manifold with or without boundary, and let $\boldsymbol{\omega}$ be a principal connection on $\mathbf{P}\to M$. A section
$\mathbf{u}\in\Gamma(\operatorname{Ad}(\mathbf{P}))$ is
\textbf{$\boldsymbol{\omega}$-parallel} if its equivariant map
$\widehat u\colon \mathbf{P}\to G$ is constant along every horizontal curve.
Denote the group of these sections by
$\Gamma^{\boldsymbol{\omega}}_{\mathrm{par}}
(\operatorname{Ad}(\mathbf{P}))$.
\end{definition}

\begin{theorem}[Parallel sections and the centralizer of holonomy]
\label{teo:secciones-paralelas-centralizador-holonomia}
\index{centralizer!of holonomy}
Let $M$ be a connected smooth manifold with or without boundary, let
$\boldsymbol{\omega}$ be a principal connection on $\mathbf{P}\to M$, and fix
$p\in \mathbf{P}$. Evaluation
\[
 \operatorname{ev}_p\colon
 \Gamma^{\boldsymbol{\omega}}_{\mathrm{par}}
 (\operatorname{Ad}(\mathbf{P}))
 \longrightarrow G,\qquad \mathbf{u}\longmapsto \widehat u(p),
\]
is a group isomorphism onto
$Z_G(\operatorname{Hol}_p(\boldsymbol{\omega}))$.
\end{theorem}
\begin{proof}
Let $\mathbf{u}$ be parallel, let $\widehat u$ be its equivariant representative, and set $c=\widehat u(p)$. If a loop has a horizontal lift from $p$ to $pg$, constancy along horizontal curves and equivariance give
\[
 c=\widehat u(p)=\widehat u(pg)
 =g^{-1}\widehat u(p)g=g^{-1}cg.
\]
Therefore $c$ centralizes
$\operatorname{Hol}_p(\boldsymbol{\omega})$.

Conversely, let
$c\in Z_G(\operatorname{Hol}_p(\boldsymbol{\omega}))$. For
$q\in \mathbf{P}_y$, choose a curve $\gamma$ from $\pi(p)$ to $y$ and write uniquely
\[
 q=\mathcal P^{\boldsymbol{\omega}}_\gamma(p)h.
\]
Define $\widehat u(q):=h^{-1}ch$. If another curve produces
$q=\mathcal P^{\boldsymbol{\omega}}_{\gamma'}(p)h'$, then
$\mathcal P^{\boldsymbol{\omega}}_{\gamma'}(p)
=\mathcal P^{\boldsymbol{\omega}}_\gamma(p)a$ for some
$a\in\operatorname{Hol}_p(\boldsymbol{\omega})$ and $h'=a^{-1}h$. Since $ac=ca$, the two formulas for $\widehat u(q)$ agree. The construction satisfies
$\widehat u(qk)=k^{-1}\widehat u(q)k$ and is constant along horizontal lifts. In a trivializing neighborhood, choose paths depending smoothly on the point; smooth dependence of transport shows that $\widehat u$ is smooth.

The construction is inverse to evaluation. It also shows that a parallel section equal to $e$ at $p$ is the identity section, so evaluation is injective. Since multiplication is pointwise, it is a group isomorphism.
\end{proof}

\begin{definition}[Gauge group and right action on connections]
\label{def:grupo-gauge-accion-derecha-conexiones}
Let $M$ be a smooth manifold with or without boundary, and let
$\mathbf{P}\to M$ be a right principal bundle. A
\textit{gauge transformation} $\mathbf{u}$ is a section of
$\operatorname{Ad}(\mathbf{P})$ or, equivalently, is represented by a smooth map $\widehat u\colon \mathbf{P}\to G$ satisfying
\[
 \widehat u(p\cdot g)=g^{-1}\widehat u(p)g.
\]
The corresponding principal automorphism is
\[
 \Phi_{\mathbf{u}}\colon \mathbf{P}\longrightarrow \mathbf{P},\qquad
 \Phi_{\mathbf{u}}(p):=p\cdot\widehat u(p).
\]
Define the right gauge action on connections by
\[
 \boldsymbol{\omega}\mathbin{\triangleleft}\mathbf{u}
 :=\Phi_{\mathbf{u}}^*\boldsymbol{\omega}.
\]
\end{definition}

\begin{proposition}[Three models of the gauge group]
\label{prop:tres-modelos-grupo-gauge}
Let $M$ be a smooth manifold with or without boundary, and let
$\mathbf{P}\to M$ be a right principal bundle.
There is a group correspondence between:
\begin{enumerate}[label=(\alph*)]
\item principal automorphisms $\Phi\colon \mathbf{P}\to \mathbf{P}$ covering
$\operatorname{id}_M$;
\item smooth equivariant maps
$\widehat u\colon \mathbf{P}\to G$ with
$\widehat u(pg)=g^{-1}\widehat u(p)g$;
\item smooth sections $\mathbf{u}$ of
$\operatorname{Ad}(\mathbf{P})$.
\end{enumerate}
The correspondence between the first two models is
$\Phi_{\mathbf{u}}(p)=p\cdot\widehat u(p)$. For a connection
$\boldsymbol{\omega}$, the stabilizer is
\[
 \operatorname{Stab}(\boldsymbol{\omega})
 :=\{\mathbf{u}\in\mathcal G(\mathbf{P})\mid
 \boldsymbol{\omega}\mathbin{\triangleleft}\mathbf{u}
 =\boldsymbol{\omega}\}.
\]
If $\boldsymbol{\xi}\in\Gamma(\operatorname{ad}(\mathbf{P}))$ and
$\mathbf{u}_t$ is the gauge transformation with equivariant representative
$\widehat u_t(p)=\exp_G(t\widehat\xi(p))$, the differential at the identity of the right gauge action is
\begin{equation}
\label{eq:accion-infinitesimal-gauge-principal}
\left.\frac{d}{dt}\right|_{0}
 \bigl(\boldsymbol{\omega}\mathbin{\triangleleft}\mathbf{u}_t\bigr)
 =d_{\boldsymbol{\omega}}\boldsymbol{\xi};
\end{equation}
therefore the space of infinitesimal directions annihilating the orbit map at
$\boldsymbol{\omega}$ is $\ker d_{\boldsymbol{\omega}}$. The identification of this kernel with the Lie algebra of the Sobolev stabilizer will be proved later, once that Lie group structure has been constructed.
\end{proposition}

\begin{proof}
If $\Phi$ covers the identity, freeness and transitivity of the action on each fiber give a unique $\widehat u(p)$ with
$\Phi(p)=p\cdot\widehat u(p)$. From
$\Phi(pg)=\Phi(p)g$, we obtain
$(p\cdot g)\cdot\widehat u(p\cdot g)
=(p\cdot\widehat u(p))\cdot g$ and, by freeness,
$\widehat u(pg)=g^{-1}\widehat u(p)g$. In trivializations, the function
$\widehat u$ is smooth. The converse is checked by substituting the equivariance law, and
Proposition~\ref{prop:secciones-haz-asociado-equivariantes}, applied to the conjugation action, gives the third model.

In a local section, let $\mathbf{A}$ be the connection potential and
$\xi_U$ the local representative of $\boldsymbol{\xi}$. Differentiating
$u_t^{-1}\mathbf{A}u_t+u_t^{-1}du_t$ for
$u_t=\exp_G(t\xi_U)$ gives
$d\xi_U+[\mathbf{A},\xi_U]=d_{\mathbf{A}}\xi_U$. Thus the kernel of the differential of the orbit map at the identity is exactly
$\ker d_{\boldsymbol{\omega}}$, which is the infinitesimal assertion in the statement.
\end{proof}

\begin{proposition}[Gauge formulas]
\label{prop:formulas-gauge-conexiones-principales}
Let $M$ be a smooth manifold with or without boundary, let $\boldsymbol{\omega}$ be a principal connection on $\mathbf{P}\to M$, and let
$\mathbf{u},\mathbf{v}\in\mathcal G(\mathbf{P})$. With pointwise multiplication of gauge transformations,
\[
 \Phi_{\mathbf{u}}\circ\Phi_{\mathbf{v}}=\Phi_{\mathbf{u}\mathbf{v}},
 \qquad
 (\boldsymbol{\omega}\mathbin{\triangleleft}\mathbf{u})
 \mathbin{\triangleleft}\mathbf{v}
 =
 \boldsymbol{\omega}\mathbin{\triangleleft}(\mathbf{u}\mathbf{v}).
\]
If $\mathbf{s}\colon U\to \mathbf{P}$ is a local section,
$\mathbf{A}=\mathbf{s}^*\boldsymbol{\omega}$, and
$u_U=\widehat u\circ\mathbf{s}$, then
\begin{align}
 \mathbf{A}\mathbin{\triangleleft}u_U
 &:=
 \mathbf{s}^*(\boldsymbol{\omega}\mathbin{\triangleleft}\mathbf{u})
 =
 \operatorname{Ad}_{u_U^{-1}}\mathbf{A}+u_U^{-1}du_U,
 \label{eq:accion-gauge-derecha-potencial}\\
 \mathbf{F}_{\mathbf{A}\mathbin{\triangleleft}u_U}
 &=
 \operatorname{Ad}_{u_U^{-1}}\mathbf{F}_{\mathbf{A}}.
 \label{eq:accion-gauge-curvatura}
\end{align}
For a matrix group,
\[
 \mathbf{A}\mathbin{\triangleleft}u
 =u^{-1}\mathbf{A}u+u^{-1}du,
 \qquad
 \mathbf{F}_{\mathbf{A}\mathbin{\triangleleft}u}
 =u^{-1}\mathbf{F}_{\mathbf{A}}u.
\]
\end{proposition}

\begin{proof}
Equivariance of the representatives $\widehat u$ and $\widehat v$ gives
\[
 \Phi_{\mathbf{u}}(\Phi_{\mathbf{v}}(p))
 =
 p\,\widehat v(p)\,\widehat u(p\widehat v(p))
 =
 p\,\widehat v(p)\,\widehat v(p)^{-1}
 \widehat u(p)\widehat v(p)
 =
 p\,\widehat u(p)\widehat v(p)
 =
 \Phi_{\mathbf{u}\mathbf{v}}(p).
\]
By contravariance of pullback,
\[
 \Phi_{\mathbf{v}}^*\Phi_{\mathbf{u}}^*\boldsymbol{\omega}
 =
 (\Phi_{\mathbf{u}}\circ\Phi_{\mathbf{v}})^*\boldsymbol{\omega}
 =
 \Phi_{\mathbf{u}\mathbf{v}}^*\boldsymbol{\omega},
\]
which proves the right action law.

Moreover,
\[
 \Phi_{\mathbf{u}}\circ\mathbf{s}=\mathbf{s}\,u_U.
\]
The local transformation formula
\eqref{eq:transformacion-potencial-conexion-principal}, applied to the sections $\mathbf{s}$ and $\mathbf{s}\,u_U$, gives
\eqref{eq:accion-gauge-derecha-potencial}. Horizontality and equivariance of curvature, or equivalently
\eqref{eq:transformacion-curvatura-principal}, give
\eqref{eq:accion-gauge-curvatura}.
\end{proof}

\begin{theorem}[Smooth stabilizer and holonomy]
\label{teo:estabilizador-suave-centralizador-holonomia}
\index{gauge stabilizer!and holonomy}
Let $M$ be a connected smooth manifold with or without boundary, let
$\boldsymbol{\omega}$ be a principal connection on $\mathbf{P}\to M$, and fix $p\in \mathbf{P}$. For a smooth gauge transformation
$\mathbf{u}$ with equivariant representative $\widehat u$, the following are equivalent:
\begin{enumerate}[label=(\alph*)]
\item
$\boldsymbol{\omega}\mathbin{\triangleleft}\mathbf{u}
=\boldsymbol{\omega}$;
\item the equivariant map $\widehat u\colon \mathbf{P}\to G$ is constant along every horizontal curve;
\item the section $\mathbf{u}$ of $\operatorname{Ad}(\mathbf{P})$ is
$\boldsymbol{\omega}$-parallel.
\end{enumerate}
Consequently, evaluation at $p$ induces an isomorphism
\begin{equation}
 \operatorname{Stab}_{C^\infty}(\boldsymbol{\omega})
 \xrightarrow{\ \cong\ }
 Z_G(\operatorname{Hol}_p(\boldsymbol{\omega})).
 \label{eq:estabilizador-suave-centralizador-holonomia}
\end{equation}
Moreover, evaluation identifies the infinitesimal directions:
\begin{equation}
 \ker\bigl(d_{\boldsymbol{\omega}}\colon
 \Omega^0(M;\operatorname{ad}(\mathbf{P}))
 \to\Omega^1(M;\operatorname{ad}(\mathbf{P}))\bigr)
 \cong
 \operatorname{Lie}
 Z_G(\operatorname{Hol}_p(\boldsymbol{\omega})).
 \label{eq:nucleo-dA-algebra-centralizador-holonomia}
\end{equation}
\end{theorem}
\begin{proof}
On a horizontal vector $X\in \mathbf{H}_q\mathbf{P}$, the global transformation formula for a connection reduces to
\[
 (\Phi_{\mathbf{u}}^*\boldsymbol{\omega})_q(X)
 =\theta^L_{\widehat u(q)}(d\widehat u_qX),
\]
because $\boldsymbol{\omega}_q(X)=0$. Therefore
$\Phi_{\mathbf{u}}^*\boldsymbol{\omega}=\boldsymbol{\omega}$ implies
$d\widehat u(X)=0$ for every horizontal vector.
Conversely, this last condition makes both connection forms vanish on $\mathbf{H}$; on vertical vectors, both reproduce the infinitesimal generator. Since
$T\mathbf{P}=\mathbf{H}\oplus\mathbf{V}$, the forms agree. This proves the equivalence of (a)--(c).

The isomorphism \eqref{eq:estabilizador-suave-centralizador-holonomia} is now
Theorem~\ref{teo:secciones-paralelas-centralizador-holonomia}.
For the infinitesimal statement, the formula for the induced connection on
$\operatorname{ad}(\mathbf{P})$ shows that
$d_{\boldsymbol{\omega}}\boldsymbol{\xi}=0$ if and only if the equivariant map $\widehat\xi\colon \mathbf{P}\to\mathfrak g$ is constant along horizontal curves. The same evaluation argument identifies it with the vectors fixed by
$\operatorname{Ad}(\operatorname{Hol}_p(\boldsymbol{\omega}))$, which form the Lie algebra of the centralizer by
\eqref{eq:algebra-lie-centralizador-fijos-ad}.
\end{proof}

\section{Covariant differentiation along curves and parallel transport}

In many geometric constructions, differentiating sections defined on open sets of the manifold is not enough. We often need to study objects defined only along a curve. This occurs, for example, when comparing elements belonging to different fibers, studying the variation of a section along a path, or, in the particular case of the tangent bundle, defining the covariant acceleration of a curve.

The difficulty is that, if $\pi\colon \mathbf{E}\longrightarrow M$ is a smooth vector bundle and $\gamma\colon I\longrightarrow M$ is a smooth curve, a section of $\mathbf{E}$ along $\gamma$ assigns to each $t\in I$ a vector in the fiber $\mathbf{E}_{\gamma(t)}$. Since these fibers are, in principle, different vector spaces, differentiating such a section by directly subtracting its values makes no sense. A connection resolves this problem: it gives an intrinsic way of differentiating sections in tangent directions and thereby induces a derivative for sections defined only along the curve.

Before formalizing this construction, we recall the pullback of a vector bundle.

\begin{definition}[Pullbacks of vector bundles]\label{def:variedades-riemannianas-pullback-de-haces-vectoriales}\index{pullbacks of vector bundles}
Let $M$ and $N$ be smooth manifolds with or without boundary, let $\pi\colon \mathbf{E}\longrightarrow M$ be a smooth vector bundle of rank $k$, and let $f\colon N\longrightarrow M$ be a smooth map. Define the \textit{pullback} of $\mathbf{E}$ by $f$ to be the set
\[
f^*\mathbf{E}:=\{(y,e)\in N\times \mathbf{E}\mid \pi(e)=f(y)\}.
\]
Equip $f^*\mathbf{E}$ with the projection $\pi_{f^*\mathbf{E}}\colon f^*\mathbf{E}\longrightarrow N$ given by $\pi_{f^*\mathbf{E}}(y,e):=y$. For each $y\in N$, the fiber of $f^*\mathbf{E}$ over $y$ is given by $(f^*\mathbf{E})_y=\{(y,e)\in f^*\mathbf{E}\mid e\in \mathbf{E}_{f(y)}\}$. Using the natural identification $(y,e)\mapsto e$, we identify $(f^*\mathbf{E})_y$ with $\mathbf{E}_{f(y)}$.
\end{definition}

\begin{proposition}\label{prop: pullback haz vectorial suave}
Let $M$ and $N$ be smooth manifolds with or without boundary, let
$\pi\colon \mathbf{E}\longrightarrow M$ be a smooth vector bundle of rank $k$,
and let $f\colon N\longrightarrow M$ be a smooth map. Then
$f^*\mathbf{E}\longrightarrow N$ admits a natural structure of a smooth vector bundle of rank $k$.
\end{proposition}

\begin{proof}
Recall that $f^*\mathbf{E}=\{(y,e)\in N\times \mathbf{E}\mid \pi(e)=f(y)\}$, with projection $\pi_{f^*\mathbf{E}}(y,e)=y$. For each $y\in N$, the fiber of $f^*\mathbf{E}$ over $y$ is given by $(f^*\mathbf{E})_y=\{(y,e)\mid e\in \mathbf{E}_{f(y)}\}$, which we identify with $\mathbf{E}_{f(y)}$ through $(y,e)\mapsto e$.

Let $\{U_\alpha\}_{\alpha\in A}$ be a trivializing cover of $\mathbf{E}$, with trivializations $\Phi_\alpha\colon \pi^{-1}(U_\alpha)\longrightarrow U_\alpha\times\mathbb K^k$. Then $\{f^{-1}(U_\alpha)\}_{\alpha\in A}$ is an open cover of $N$ because $f$ is continuous.

For each $\alpha$, define
\[
\widetilde{\Phi}_\alpha:
\pi_{f^*\mathbf{E}}^{-1}(f^{-1}(U_\alpha))
\longrightarrow
f^{-1}(U_\alpha)\times \mathbb K^k
\]
by
\[
\widetilde{\Phi}_\alpha(y,e)
:=
\bigl(y,\operatorname{pr}_2(\Phi_\alpha(e))\bigr),
\]
where $\operatorname{pr}_2\colon U_\alpha\times\mathbb K^k\longrightarrow\mathbb K^k$ is the second projection. This expression is well defined because, if $(y,e)\in f^*\mathbf{E}$ and $y\in f^{-1}(U_\alpha)$, then $f(y)\in U_\alpha$ and, since $\pi(e)=f(y)$, we have $e\in\pi^{-1}(U_\alpha)$.

Let us show that $\widetilde{\Phi}_\alpha$ is bijective. Let $(y,v)\in f^{-1}(U_\alpha)\times\mathbb K^k$. Since $f(y)\in U_\alpha$, there is a unique $e\in \mathbf{E}_{f(y)}$ such that $\Phi_\alpha(e)=(f(y),v)$. Therefore $(y,e)\in f^*\mathbf{E}$ and $\widetilde{\Phi}_\alpha(y,e)=(y,v)$. This proves surjectivity.

If $\widetilde{\Phi}_\alpha(y,e)=\widetilde{\Phi}_\alpha(y',e')$, then $y=y'$ and $\operatorname{pr}_2(\Phi_\alpha(e))=\operatorname{pr}_2(\Phi_\alpha(e'))$. Moreover, since $(y,e)$ and $(y,e')$ belong to $f^*\mathbf{E}$, we have $\pi(e)=f(y)=\pi(e')$. Therefore the first components of $\Phi_\alpha(e)$ and $\Phi_\alpha(e')$ also agree. Thus $\Phi_\alpha(e)=\Phi_\alpha(e')$. Since $\Phi_\alpha$ is injective, we obtain $e=e'$. Hence $\widetilde{\Phi}_\alpha$ is injective. In fact, its inverse is given by
\[
\widetilde{\Phi}_\alpha^{-1}(y,v)
=
\bigl(y,\Phi_\alpha^{-1}(f(y),v)\bigr).
\]

Moreover, the restriction of $\widetilde{\Phi}_\alpha$ to each fiber $(f^*\mathbf{E})_y$ is a linear isomorphism onto $\{y\}\times\mathbb K^k$, since it agrees with the trivialization of $\mathbf{E}$ on the fiber $\mathbf{E}_{f(y)}$ under the natural identification $(f^*\mathbf{E})_y\cong \mathbf{E}_{f(y)}$.

Now consider $\alpha,\beta$ such that $f^{-1}(U_\alpha)\cap f^{-1}(U_\beta)\neq\varnothing$. On $U_\alpha\cap U_\beta$, the trivializations of $\mathbf{E}$ satisfy
\[
\Phi_\alpha\circ\Phi_\beta^{-1}(p,v)
=
(p,\tau_{\alpha\beta}(p)v),
\]
where $\tau_{\alpha\beta}\colon U_\alpha\cap U_\beta\longrightarrow GL(k,\mathbb K)$ is smooth. Then, for $y\in f^{-1}(U_\alpha\cap U_\beta)$, we have
\[
\widetilde{\Phi}_\alpha\circ\widetilde{\Phi}_\beta^{-1}(y,v)
=
\bigl(y,\tau_{\alpha\beta}(f(y))v\bigr).
\]
Thus the transition functions of the system $\{\widetilde{\Phi}_\alpha\}$ are given by $\widetilde{\tau}_{\alpha\beta}(y)=\tau_{\alpha\beta}(f(y))$. Since $\tau_{\alpha\beta}$ and $f$ are smooth, $\widetilde{\tau}_{\alpha\beta}$ is also smooth.

By Lemma~\ref{lema del haz vectorial suave}, there exist a unique topology and smooth structure on $f^*\mathbf{E}$ making it a smooth vector bundle of rank $k$ over $N$, with projection $\pi_{f^*\mathbf{E}}$ and local trivializations $\widetilde{\Phi}_\alpha$.
\end{proof}

\begin{definition}[Sections along a curve]\label{def:variedades-riemannianas-secciones-a-lo-largo-de-una-curva}\index{sections along a curve}
Let $M$ be a smooth manifold with or without boundary, let
$\pi\colon \mathbf{E}\longrightarrow M$ be a smooth vector bundle, and let
$\gamma\colon I\longrightarrow M$ be a smooth curve. A section of
$\mathbf{E}$ along $\gamma$ is a section of the pullback bundle
$\gamma^*\mathbf{E}\longrightarrow I$.

Using the natural identification $(\gamma^*\mathbf{E})_t\cong \mathbf{E}_{\gamma(t)}$, we regard a section of $\mathbf{E}$ along $\gamma$ as a map $\mathbf{u}\colon I\longrightarrow \mathbf{E}$ such that $\mathbf{u}(t)\in \mathbf{E}_{\gamma(t)}$ for every $t\in I$. Under this identification, we write $\Gamma^k(\gamma^*\mathbf{E})$ for the space of sections of class $C^k$ of $\mathbf{E}$ along $\gamma$, and we write $\Gamma(\gamma^*\mathbf{E})=\Gamma^\infty(\gamma^*\mathbf{E})$.

We say that $\mathbf{u}\in\Gamma^k(\gamma^*\mathbf{E})$ is extendible if there exists a section $\widetilde{\mathbf{u}}\in\Gamma^k(\mathbf{E})$, defined on a neighborhood of $\gamma(I)$, such that $\mathbf{u}(t)=\widetilde{\mathbf{u}}(\gamma(t))$ for every $t\in I$.
\end{definition}

\begin{remark}\label{obs:variedades-riemannianas-cuando-elementos-precisamente-campos-vectoriales-clase}
When $\mathbf{E}=TM$, the elements of $\Gamma^k(\gamma^*TM)$ are precisely vector fields of class $C^k$ along $\gamma$. More generally, if $\mathbf{E}=T^{(r,s)}(TM)$ carries the connection induced by a connection on $TM$, we obtain the notion of a tensor field of type $(r,s)$ along $\gamma$.
\end{remark}

The following result introduces covariant differentiation along a curve. The essential condition is that, if a section along $\gamma$ is obtained by restricting a section $\widetilde{\mathbf{u}}$ of $\mathbf{E}$, its derivative along the curve must agree with $\nabla^{\mathbf{E}}_{\dot\gamma}\widetilde{\mathbf{u}}$. The result states that this condition uniquely determines a derivative for every section of class $C^1$ of $\gamma^*\mathbf{E}$.

\begin{theorem}[Covariant derivative along a curve]\label{teo: derivada covariante a lo largo de una curva}\index{covariant derivative along a curve}
Let $M$ be a smooth manifold with or without boundary, let
$\pi\colon \mathbf{E}\longrightarrow M$ be a smooth vector bundle with connection
$\nabla^{\mathbf{E}}$, and let $\gamma\colon I\longrightarrow M$ be a smooth curve. There exists a unique operator
\[
\displaystyle\frac{D}{dt}\colon \Gamma^1(\gamma^*\mathbf{E})\longrightarrow \Gamma^0(\gamma^*\mathbf{E})
\]
satisfying the following properties:

\begin{enumerate}[label=(\alph*)]
 \item $\displaystyle \frac{D}{dt}(a\mathbf{u}+b\mathbf{v})=a\frac{D\mathbf{u}}{dt}+b\frac{D\mathbf{v}}{dt}$ for every $a,b\in\mathbb K$ and every $\mathbf{u},\mathbf{v}\in\Gamma^1(\gamma^*\mathbf{E})$.

 \item $\displaystyle \frac{D}{dt}(f\mathbf{u})=\frac{df}{dt}\mathbf{u}+f\frac{D\mathbf{u}}{dt}$ for every $f\in C^1(I)$ and every $\mathbf{u}\in\Gamma^1(\gamma^*\mathbf{E})$.

 \item If $\mathbf{u}$ is extendible and $\widetilde{\mathbf{u}}$ is an extension of class $C^1$ of $\mathbf{u}$ to a neighborhood of $\gamma(I)$, then $\displaystyle \frac{D\mathbf{u}}{dt}(t)=\nabla^{\mathbf{E}}_{\dot\gamma(t)}\widetilde{\mathbf{u}}$ for every $t\in I$.
\end{enumerate}

If $\mathbf{u}\in\Gamma^k(\gamma^*\mathbf{E})$ with $k\geq 1$, then $\displaystyle \frac{D\mathbf{u}}{dt}\in\Gamma^{k-1}(\gamma^*\mathbf{E})$. In particular, this operator restricts to an operator $\displaystyle \frac{D}{dt}\colon \Gamma(\gamma^*\mathbf{E})\longrightarrow\Gamma(\gamma^*\mathbf{E})$.
\end{theorem}

\begin{proof}
We first construct the operator in a local trivialization and then prove independence of the choices.

Let $t_0\in I$. Choose a smooth chart $(U,\phi)$ on $M$, with local coordinates $(x^1,\dots,x^n)$, such that $\gamma(t_0)\in U$, and choose a smooth local frame $(\mathbf{e}_1,\dots,\mathbf{e}_r)$ of $\mathbf{E}$ over $U$. Since $\gamma^{-1}(U)$ is open in $I$ and contains $t_0$, there is an interval $J\subseteq I$, open in $I$, such that $t_0\in J$ and $\gamma(J)\subseteq U$.

On $J$, every section $\mathbf{u}\in\Gamma^1(\gamma^*\mathbf{E})$ can be written uniquely as $\mathbf{u}(t)=u^\alpha(t)\mathbf{e}_\alpha(\gamma(t))$, where $u^1,\dots,u^r$ are functions of class $C^1$ on $J$. Define the local connection coefficients by $\nabla^{\mathbf{E}}_{\boldsymbol{\partial}_i}\mathbf{e}_\beta=\Gamma^\alpha_{i\beta}\mathbf{e}_\alpha$. Then define locally
\begin{equation}\label{eq: derivada covariante curva haz local}
\displaystyle\frac{D\mathbf{u}}{dt}(t)
=
\left(
\displaystyle\frac{du^\alpha}{dt}(t)
+
\displaystyle\frac{d\gamma^i}{dt}(t)\Gamma^\alpha_{i\beta}(\gamma(t))u^\beta(t)
\right)\mathbf{e}_\alpha(\gamma(t)).
\end{equation}

We show that this definition is independent of the local frame. Let $(\widetilde{\mathbf{e}}_1,\dots,\widetilde{\mathbf{e}}_r)$ be another local frame of $\mathbf{E}$ over the same open set, and write $\widetilde{\mathbf{e}}_a=B^\beta_a \mathbf{e}_\beta$, where $B=(B^\beta_a)$ is a matrix of smooth functions with values in $GL(r,\mathbb K)$. If
\[
\mathbf{u}(t)=u^\alpha(t)\mathbf{e}_\alpha(\gamma(t))
=
\widetilde u^a(t)\widetilde{\mathbf{e}}_a(\gamma(t)),
\]
then $u^\alpha(t)=B^\alpha_a(\gamma(t))\widetilde u^a(t)$. On the other hand, if $\nabla^{\mathbf{E}}_{\boldsymbol{\partial}_i}\widetilde{\mathbf{e}}_b=\widetilde\Gamma^a_{ib}\widetilde{\mathbf{e}}_a$, then, using $\widetilde{\mathbf{e}}_b=B^\beta_b \mathbf{e}_\beta$, we obtain
\[
\nabla^{\mathbf{E}}_{\boldsymbol{\partial}_i}\widetilde{\mathbf{e}}_b
=
(\partial_iB^\alpha_b)\mathbf{e}_\alpha
+
B^\beta_b\nabla^{\mathbf{E}}_{\boldsymbol{\partial}_i}\mathbf{e}_\beta
=
(\partial_iB^\alpha_b+B^\beta_b\Gamma^\alpha_{i\beta})\mathbf{e}_\alpha.
\]
Since also
\[
\nabla^{\mathbf{E}}_{\boldsymbol{\partial}_i}\widetilde{\mathbf{e}}_b
=
\widetilde\Gamma^a_{ib}\widetilde{\mathbf{e}}_a
=
\widetilde\Gamma^a_{ib}B^\alpha_a \mathbf{e}_\alpha,
\]
we have
\[
\partial_iB^\alpha_b+B^\beta_b\Gamma^\alpha_{i\beta}
=
B^\alpha_a\widetilde\Gamma^a_{ib}.
\]
Differentiating $u^\alpha=B^\alpha_a(\gamma)\widetilde u^a$ with respect to $t$ gives
\[
\displaystyle\frac{du^\alpha}{dt}
=
\displaystyle\frac{d\gamma^i}{dt}(\partial_iB^\alpha_a)(\gamma)\widetilde u^a
+
B^\alpha_a(\gamma)\displaystyle\frac{d\widetilde u^a}{dt}.
\]
Therefore,
\[
\displaystyle\frac{du^\alpha}{dt}
+
\displaystyle\frac{d\gamma^i}{dt}\Gamma^\alpha_{i\beta}(\gamma)u^\beta
=
B^\alpha_a(\gamma)\displaystyle\frac{d\widetilde u^a}{dt}
+
\displaystyle\frac{d\gamma^i}{dt}
\bigl(\partial_iB^\alpha_b+B^\beta_b\Gamma^\alpha_{i\beta}\bigr)(\gamma)\widetilde u^b.
\]
Using the preceding relation between the connection coefficients yields
\[
\displaystyle\frac{du^\alpha}{dt}
+
\displaystyle\frac{d\gamma^i}{dt}\Gamma^\alpha_{i\beta}(\gamma)u^\beta
=
B^\alpha_a(\gamma)
\left(
\displaystyle\frac{d\widetilde u^a}{dt}
+
\displaystyle\frac{d\gamma^i}{dt}\widetilde\Gamma^a_{ib}(\gamma)\widetilde u^b
\right).
\]
Multiplying by $\mathbf{e}_\alpha(\gamma(t))$ and using $B^\alpha_a \mathbf{e}_\alpha=\widetilde{\mathbf{e}}_a$ gives exactly the expression for \eqref{eq: derivada covariante curva haz local} in the frame $\widetilde{\mathbf{e}}_a$. Thus the definition is independent of the local frame.

Since the resulting expression is intrinsic with respect to the frame, it also agrees on overlaps of local trivializations. It therefore defines a global section $\displaystyle \frac{D\mathbf{u}}{dt}\in\Gamma^0(\gamma^*\mathbf{E})$. If $\mathbf{u}\in\Gamma^k(\gamma^*\mathbf{E})$, the local formula shows that $\displaystyle \frac{D\mathbf{u}}{dt}\in\Gamma^{k-1}(\gamma^*\mathbf{E})$.

Linearity over $\mathbb K$ is checked directly in the local formula. For the Leibniz rule, let $f\in C^1(I)$. The components of $f\mathbf{u}$ are $fu^\alpha$, and then
\[
\displaystyle\frac{D(f\mathbf{u})}{dt}
=
\left(
\displaystyle\frac{d(fu^\alpha)}{dt}
+
\displaystyle\frac{d\gamma^i}{dt}\Gamma^\alpha_{i\beta}(\gamma)fu^\beta
\right)\mathbf{e}_\alpha(\gamma(t)).
\]
Since $\displaystyle \frac{d(fu^\alpha)}{dt}=\frac{df}{dt}u^\alpha+f\frac{du^\alpha}{dt}$, we obtain $\displaystyle \frac{D(f\mathbf{u})}{dt}=\frac{df}{dt}\mathbf{u}+f\frac{D\mathbf{u}}{dt}$.

We now prove compatibility with extensions. Suppose that $\mathbf{u}(t)=\widetilde{\mathbf{u}}(\gamma(t))$, where $\widetilde{\mathbf{u}}$ is a section of class $C^1$ of $\mathbf{E}$ defined on a neighborhood of $\gamma(I)$. In a local frame, write $\widetilde{\mathbf{u}}=\widetilde u^\alpha \mathbf{e}_\alpha$. The components of $\mathbf{u}$ along $\gamma$ are then $\widetilde u^\alpha(\gamma(t))$. By the chain rule,
\[
\displaystyle\frac{d}{dt}\bigl(\widetilde u^\alpha(\gamma(t))\bigr)
=
\displaystyle\frac{d\gamma^i}{dt}(t)\partial_i\widetilde u^\alpha(\gamma(t)).
\]
Substituting into \eqref{eq: derivada covariante curva haz local} gives
\[
\displaystyle\frac{D\mathbf{u}}{dt}(t)
=
\displaystyle\frac{d\gamma^i}{dt}(t)
\left(
\partial_i\widetilde u^\alpha(\gamma(t))
+
\Gamma^\alpha_{i\beta}(\gamma(t))\widetilde u^\beta(\gamma(t))
\right)\mathbf{e}_\alpha(\gamma(t)).
\]
The expression in parentheses is the local expression for $\nabla^{\mathbf{E}}_{\boldsymbol{\partial}_i}\widetilde{\mathbf{u}}$ in the frame $\mathbf{e}_\alpha$. Therefore $\displaystyle \frac{D\mathbf{u}}{dt}(t)=\nabla^{\mathbf{E}}_{\dot\gamma(t)}\widetilde{\mathbf{u}}$.

Finally, we prove uniqueness. Suppose that $\mathcal D$ is another operator satisfying the three properties in the statement. On an interval $J$ with $\gamma(J)\subseteq U$, write $\mathbf{u}(t)=u^\alpha(t)\mathbf{e}_\alpha(\gamma(t))$. Since each $\mathbf{e}_\alpha\circ\gamma$ is extendible, the properties of $\mathcal D$ imply
\[
\mathcal D \mathbf{u}
=
\displaystyle\frac{du^\alpha}{dt}\mathbf{e}_\alpha(\gamma(t))
+
u^\alpha(t)\mathcal D(\mathbf{e}_\alpha\circ\gamma).
\]
By compatibility with extensions, $\mathcal D(\mathbf{e}_\alpha\circ\gamma)(t)=\nabla^{\mathbf{E}}_{\dot\gamma(t)}\mathbf{e}_\alpha$. Using $\dot\gamma(t)=\frac{d\gamma^i}{dt}(t)\boldsymbol{\partial}_i|_{\gamma(t)}$ and $\nabla^{\mathbf{E}}_{\boldsymbol{\partial}_i}\mathbf{e}_\beta=\Gamma^\alpha_{i\beta}\mathbf{e}_\alpha$ gives precisely formula \eqref{eq: derivada covariante curva haz local}. Therefore $\mathcal D=\frac{D}{dt}$ locally. Since $t_0$ was arbitrary, $\mathcal D=\frac{D}{dt}$ on all of $I$.
\end{proof}

\begin{remark}\label{obs:variedades-riemannianas-consideran-conexiones-inducidas-productos-tensoriales-partic}
For the induced connections on $\mathbf{E}^*$, on tensor products, and in particular on tensor bundles constructed from $TM$, covariant differentiation along curves commutes with contractions and satisfies the Leibniz rule for the tensor product. In particular, taking $\mathbf{E}=T^{(r,s)}(TM)$ gives covariant differentiation of tensor fields of type $(r,s)$ along a curve.
\end{remark}

\begin{proposition}[Local expression for the covariant derivative along a curve]\label{prop: expresion local derivada covariante curvas}\index{derivative!covariant along a curve!local expression}
Let $M$ be a smooth manifold with or without boundary, let
$\pi\colon \mathbf{E}\longrightarrow M$ be a smooth vector bundle of rank $r$ with connection $\nabla^{\mathbf{E}}$, and let $\gamma\colon I\longrightarrow M$ be a smooth curve. Suppose that $\gamma(J)\subseteq U$, where $(U,\phi)$ is a smooth chart on $M$ with coordinates $(x^1,\dots,x^n)$, and let
$(\mathbf{e}_1,\dots,\mathbf{e}_r)$ be a smooth local frame of $\mathbf{E}$ over $U$.

If $\mathbf{u}(t)=u^\alpha(t)\mathbf{e}_\alpha(\gamma(t))$, with $t\in J$, and the local connection coefficients are given by $\nabla^{\mathbf{E}}_{\boldsymbol{\partial}_i}\mathbf{e}_\beta=\Gamma^\alpha_{i\beta}\mathbf{e}_\alpha$, then
\[
\displaystyle\frac{D\mathbf{u}}{dt}(t)
=
\left(
\displaystyle\frac{du^\alpha}{dt}(t)
+
\displaystyle\frac{d\gamma^i}{dt}(t)\Gamma^\alpha_{i\beta}(\gamma(t))u^\beta(t)
\right)\mathbf{e}_\alpha(\gamma(t)).
\]
\end{proposition}

\begin{proof}
Write $\mathbf{u}=u^\alpha \mathbf{e}_\alpha\circ\gamma$. By the Leibniz rule,
\[
\displaystyle\frac{D\mathbf{u}}{dt}
=
\displaystyle\frac{du^\alpha}{dt}\mathbf{e}_\alpha(\gamma(t))
+
u^\alpha\displaystyle\frac{D(\mathbf{e}_\alpha\circ\gamma)}{dt}.
\]
Since $\mathbf{e}_\alpha\circ\gamma$ is extendible, we have $\displaystyle \frac{D(\mathbf{e}_\alpha\circ\gamma)}{dt}(t)=\nabla^{\mathbf{E}}_{\dot\gamma(t)}\mathbf{e}_\alpha$. Using $\dot\gamma(t)=\frac{d\gamma^i}{dt}(t)\boldsymbol{\partial}_i|_{\gamma(t)}$ and $\nabla^{\mathbf{E}}_{\boldsymbol{\partial}_i}\mathbf{e}_\beta=\Gamma^\alpha_{i\beta}\mathbf{e}_\alpha$, we obtain
\[
\displaystyle\frac{D(\mathbf{e}_\beta\circ\gamma)}{dt}(t)
=
\displaystyle\frac{d\gamma^i}{dt}(t)\Gamma^\alpha_{i\beta}(\gamma(t))\mathbf{e}_\alpha(\gamma(t)).
\]
Substituting this expression into the formula for $\frac{D\mathbf{u}}{dt}$ gives
\[
\displaystyle\frac{D\mathbf{u}}{dt}(t)
=
\displaystyle\frac{du^\alpha}{dt}(t)\mathbf{e}_\alpha(\gamma(t))
+
u^\beta(t)\displaystyle\frac{d\gamma^i}{dt}(t)\Gamma^\alpha_{i\beta}(\gamma(t))\mathbf{e}_\alpha(\gamma(t)).
\]
Grouping the coefficients of $\mathbf{e}_\alpha(\gamma(t))$ yields the desired formula.
\end{proof}

The preceding result shows that the covariant derivative of a section along a curve is calculated using a local system of ordinary functions. In particular, if a global section is restricted to a curve, its derivative along the curve depends only on the values of the section along that curve and on the curve's velocity.

\begin{proposition}\label{prop: derivada depende solo de la curva}
Let $M$ be a smooth manifold with or without boundary, and let
$\pi\colon \mathbf{E}\longrightarrow M$ be a smooth vector bundle with connection
$\nabla^{\mathbf{E}}$. Let $p\in M$, $v\in T_pM$, and let
$\mathbf{s},\widetilde{\mathbf{s}}$ be sections of class $C^1$ of $\mathbf{E}$ defined on a neighborhood of $p$. Suppose that there exist a smooth curve
$\gamma\colon I\longrightarrow M$ and a point $t_0\in I$ such that
$\gamma(t_0)=p$, $\dot\gamma(t_0)=v$, and
$\mathbf{s}(\gamma(t))=\widetilde{\mathbf{s}}(\gamma(t))$ for every $t$ in a neighborhood of $t_0$. Then
$\nabla^{\mathbf{E}}_v\mathbf{s}=\nabla^{\mathbf{E}}_v\widetilde{\mathbf{s}}$.
\end{proposition}

\begin{proof}
Define a section $\mathbf{u}$ of $\gamma^*\mathbf{E}$ near $t_0$ by
$\mathbf{u}(t)=\mathbf{s}(\gamma(t))=\widetilde{\mathbf{s}}(\gamma(t))$.
Then $\mathbf{s}$ and $\widetilde{\mathbf{s}}$ are two local extensions of $\mathbf{u}$. By compatibility with extensions in
Theorem~\ref{teo: derivada covariante a lo largo de una curva}, we have
\[
\nabla^{\mathbf{E}}_{\dot\gamma(t_0)}\mathbf{s}
=
\displaystyle\frac{D\mathbf{u}}{dt}(t_0)
=
\nabla^{\mathbf{E}}_{\dot\gamma(t_0)}\widetilde{\mathbf{s}}.
\]
Since $\dot\gamma(t_0)=v$, we conclude that
$\nabla^{\mathbf{E}}_v\mathbf{s}=\nabla^{\mathbf{E}}_v\widetilde{\mathbf{s}}$.
\end{proof}

\begin{definition}[Parallel sections along a curve]\label{def:variedades-riemannianas-secciones-paralelas-a-lo-largo-de-una-curva}\index{parallel sections along a curve}
Let $M$ be a smooth manifold with or without boundary, let
$\pi\colon \mathbf{E}\longrightarrow M$ be a smooth vector bundle with connection
$\nabla^{\mathbf{E}}$, and let $\gamma\colon I\longrightarrow M$ be a smooth curve. A section
$\mathbf{u}\in\Gamma^1(\gamma^*\mathbf{E})$ is said to be parallel along
$\gamma$ if $\displaystyle \frac{D\mathbf{u}}{dt}\equiv0$.
\end{definition}

In a local frame, the parallelism condition becomes a linear system of ordinary differential equations. If
$\mathbf{u}(t)=u^\alpha(t)\mathbf{e}_\alpha(\gamma(t))$, then
$\mathbf{u}$ is parallel along $\gamma$ if and only if
\begin{equation}\label{eq: transporte paralelo haz}
\displaystyle\frac{du^\alpha}{dt}(t)
=
-\displaystyle\frac{d\gamma^i}{dt}(t)\Gamma^\alpha_{i\beta}(\gamma(t))u^\beta(t),
\qquad \alpha\in\{1,\dots,r\}.
\end{equation}
This observation reduces existence and uniqueness of parallel sections to a classical existence and uniqueness theorem for linear systems of ordinary differential equations.

\begin{proposition}\label{prop: seccion paralela iff restricciones paralelas}
Let $M$ be a smooth manifold with or without boundary, let
$\pi\colon \mathbf{E}\longrightarrow M$ be a smooth vector bundle with connection
$\nabla^{\mathbf{E}}$, and let $\mathbf{s}\in\Gamma^1(\mathbf{E})$. Then
$\nabla^{\mathbf{E}} \mathbf{s}=0$ if and only if, for every smooth curve
$\gamma\colon I\longrightarrow M$, the section $\mathbf{s}\circ\gamma$ of
$\gamma^*\mathbf{E}$ is parallel along $\gamma$.
\end{proposition}

\begin{proof}
First suppose that $\nabla^{\mathbf{E}} \mathbf{s}=0$. For every smooth curve $\gamma\colon I\longrightarrow M$, the section $\mathbf{s}\circ\gamma$ is extendible by $\mathbf{s}$. By Theorem~\ref{teo: derivada covariante a lo largo de una curva}, we have $\displaystyle \frac{D(\mathbf{s}\circ\gamma)}{dt}(t)=\nabla^{\mathbf{E}}_{\dot\gamma(t)}\mathbf{s}=0$. Thus $\mathbf{s}\circ\gamma$ is parallel along $\gamma$.

Conversely, suppose that $\mathbf{s}\circ\gamma$ is parallel along every smooth curve $\gamma$. Let $p\in M$ and $v\in T_pM$. Choose a smooth curve $\gamma$ defined near $0$ such that $\gamma(0)=p$ and $\dot\gamma(0)=v$. Then
\[
\nabla^{\mathbf{E}}_v \mathbf{s}
=
\nabla^{\mathbf{E}}_{\dot\gamma(0)}\mathbf{s}
=
\displaystyle\frac{D(\mathbf{s}\circ\gamma)}{dt}(0)
=
0.
\]
Since $p$ and $v$ were arbitrary, we conclude that $\nabla^{\mathbf{E}} \mathbf{s}=0$.
\end{proof}

\begin{theorem}[Existence, uniqueness, and smoothness for linear ODEs]\label{teo: existencia y unicidad para ODEs lineales}\index{existence uniqueness and smoothness for linear ODEs@existence, uniqueness, and smoothness for linear ODEs}
Let $I\subseteq\mathbb R$ be an open interval and, for $1\leq\beta,\alpha\leq r$, let $A^\alpha_\beta\colon I\longrightarrow\mathbb K$ be a smooth function. For every $t_0\in I$ and every initial vector $c=(c^1,\dots,c^r)\in\mathbb K^r$, the linear initial value problem
\begin{equation}\label{eq: sistema de EDO lineal}
\begin{cases}
\displaystyle \frac{dU^\alpha}{dt}(t)=A^\alpha_\beta(t)U^\beta(t),\\[0.4em]
U^\alpha(t_0)=c^\alpha,
\end{cases}
\end{equation}
has a unique smooth solution defined on all of $I$. Moreover, for fixed $t_0$, the solution depends smoothly on $(t,c)\in I\times\mathbb K^r$.
\end{theorem}

\begin{proof}
We first prove the result when $t_0=0$ and $0\in I$. To simplify notation, we present the argument over $\mathbb R$. In the complex case, apply the same reasoning to the real and imaginary parts, identifying $\mathbb C^r$ with $\mathbb R^{2r}$.

Let $(x^0,x^1,\dots,x^r)$ be the standard coordinates on $I\times\mathbb R^r$. Consider the smooth vector field
\[
\mathbf{Y}
=
\boldsymbol{\partial}_{0}
+
A^\alpha_\beta(x^0)x^\beta\boldsymbol{\partial}_{\alpha}
\]
on $I\times\mathbb R^r$.

If $U(t)=(U^1(t),\dots,U^r(t))$ is a solution of \eqref{eq: sistema de EDO lineal} defined on an interval $I_0\subseteq I$, then $\eta(t)=(t,U^1(t),\dots,U^r(t))$ is an integral curve of $\mathbf{Y}$. The first component satisfies $\dot\eta^0(t)=1$, and the others satisfy
\[
\dot\eta^\alpha(t)
=
\displaystyle\frac{dU^\alpha}{dt}(t)
=
A^\alpha_\beta(t)U^\beta(t)
=
A^\alpha_\beta(\eta^0(t))\eta^\beta(t).
\]

Conversely, if $\eta(t)=(\eta^0(t),\eta^1(t),\dots,\eta^r(t))$ is an integral curve of $Y$ with $\eta(0)=(0,c^1,\dots,c^r)$, then $\dot\eta^0(t)=1$ and $\eta^0(0)=0$. Therefore $\eta^0(t)=t$. Defining $U(t)=(\eta^1(t),\dots,\eta^r(t))$ gives a solution of the system \eqref{eq: sistema de EDO lineal} with initial condition $U(0)=c$. Thus solutions of the linear system correspond to integral curves of $Y$ passing through $(0,c)$.

By Theorem~\ref{teo: fundamental de flujos}, there exists a unique maximal integral curve of $Y$ with initial condition $(0,c)$, defined on a maximal open interval $I_0\subseteq I$ containing $0$. By the preceding argument, this is equivalent to existence and uniqueness of a maximal solution $U\colon I_0\longrightarrow\mathbb R^r$ of the linear system. It remains to show that $I_0=I$.

Write $I=(a,b)$ and $I_0=(a_0,b_0)$, with $a,b,a_0,b_0\in\mathbb R\cup\{-\infty,\infty\}$. Suppose, for a contradiction, that $b_0<b$. Then $b_0$ is finite and $[0,b_0]\subseteq I$. For $t\in[0,b_0)$, we have
\[
\displaystyle\frac{d}{dt}\|U(t)\|^2
=
2\displaystyle\frac{dU^\alpha}{dt}(t)U^\alpha(t)
=
2A^\alpha_\beta(t)U^\beta(t)U^\alpha(t).
\]
Using the Cauchy--Schwarz inequality for the Frobenius norm, we obtain $\displaystyle \frac{d}{dt}\|U(t)\|^2\leq2\|A(t)\|_F\|U(t)\|^2$. Since the functions $A^\alpha_\beta$ are continuous, there exists $C>0$ such that $\|A(t)\|_F\leq C$ for every $t\in[0,b_0]$. Therefore,
\[
\displaystyle\frac{d}{dt}\left(e^{-2Ct}\|U(t)\|^2\right)
=
e^{-2Ct}
\left(
\displaystyle\frac{d}{dt}\|U(t)\|^2-2C\|U(t)\|^2
\right)
\leq 0.
\]
Thus $e^{-2Ct}\|U(t)\|^2\leq \|U(0)\|^2$ for every $t\in[0,b_0)$, and hence $\|U(t)\|\leq e^{Cb_0}\|U(0)\|$. Consequently, the corresponding integral curve remains in the compact set
\[
[0,b_0]\times \overline{B}_{\mathrm{euc}}(0,R)\subseteq I\times\mathbb R^r,
\qquad
R=e^{Cb_0}\|U(0)\|.
\]
This contradicts Lemma~\ref{lema: de escape}. Therefore $b_0=b$.

To prove that $a_0=a$, consider the reparametrized solution $W(s)=U(-s)$, which satisfies the linear system $\displaystyle \frac{dW^\alpha}{ds}(s)=-A^\alpha_\beta(-s)W^\beta(s)$. The preceding argument, applied to $W$, shows that the solution cannot cease to exist before reaching the left endpoint of $I$ either. We conclude that $I_0=I$.

Smooth dependence on $(t,c)$ follows from smooth dependence of the flow of $Y$ on time and initial conditions, since the solution $U(t)$ is the $\mathbb R^r$-component of the integral curve of $Y$ starting at $(0,c)$.

Finally, consider the general case $t_0\in I$. Define $\widetilde I=I-t_0$ and $\widetilde A^\alpha_\beta(s)=A^\alpha_\beta(s+t_0)$. The translated system
\[
\displaystyle\frac{d\widetilde U^\alpha}{ds}(s)
=
\widetilde A^\alpha_\beta(s)\widetilde U^\beta(s),
\qquad
\widetilde U^\alpha(0)=c^\alpha,
\]
has a unique smooth solution on all of $\widetilde I$. Defining $U(t)=\widetilde U(t-t_0)$ gives the solution of the original system on all of $I$. Uniqueness and smooth dependence are preserved under the change of variable $s=t-t_0$.
\end{proof}

\begin{theorem}[Existence and uniqueness of parallel transport]\label{teo: existencia y unicidad transporte paralelo}\index{existence and uniqueness of parallel transport}
Let $M$ be a smooth manifold with or without boundary, and let
$\pi\colon \mathbf{E}\longrightarrow M$ be a smooth vector bundle with connection
$\nabla^{\mathbf{E}}$. Let $\gamma\colon I\longrightarrow M$ be a smooth curve, let $t_0\in I$, and let $v\in \mathbf{E}_{\gamma(t_0)}$. Then there exists a unique parallel section
$\mathbf{u}\in\Gamma(\gamma^*\mathbf{E})$ such that $\mathbf{u}(t_0)=v$.
\end{theorem}

\begin{proof}
We first solve the problem locally. Suppose that $\gamma(J)\subseteq U$, where $J\subseteq I$ is an interval and $(\mathbf{e}_1,\dots,\mathbf{e}_r)$ is a smooth local frame of $\mathbf{E}$ over $U$. If $\mathbf{u}(t)=u^\alpha(t)\mathbf{e}_\alpha(\gamma(t))$, then, by \eqref{eq: transporte paralelo haz}, the condition $\frac{D\mathbf{u}}{dt}=0$ is equivalent to the linear system
\[
\displaystyle\frac{du^\alpha}{dt}(t)
=
-\displaystyle\frac{d\gamma^i}{dt}(t)\Gamma^\alpha_{i\beta}(\gamma(t))u^\beta(t).
\]
Defining $A^\alpha_\beta(t)=-\frac{d\gamma^i}{dt}(t)\Gamma^\alpha_{i\beta}(\gamma(t))$ gives a system of the form $\displaystyle \frac{du^\alpha}{dt}(t)=A^\alpha_\beta(t)u^\beta(t)$. Since $\gamma$ and the coefficients $\Gamma^\alpha_{i\beta}$ are smooth, the functions $A^\alpha_\beta$ are smooth on $J$. By Theorem~\ref{teo: existencia y unicidad para ODEs lineales}, for every initial condition at a point of $J$, there exists a unique smooth solution on all of $J$. Thus, on any parameter interval whose image lies in a single local trivialization of $\mathbf{E}$, there exists a unique parallel section with prescribed initial value.

We now prove global existence and uniqueness to the right of $t_0$; the left endpoint is treated by reversing the parameter. Let $b=\displaystyle\sup I$ and consider the set
\[
\mathcal S
=
\{s\in I \mid s>t_0 \text{ and there is a unique parallel section on } [t_0,s] \text{ with initial value } v\}.
\]
This set is nonempty. Choose a local trivialization of $\mathbf{E}$ on a neighborhood $U$ of $\gamma(t_0)$. Since $\gamma^{-1}(U)$ is open in $I$ and contains $t_0$, there exists $\varepsilon>0$ such that $\gamma((t_0-\varepsilon,t_0+\varepsilon)\cap I)\subseteq U$. By the local case, there is a unique parallel section on $[t_0,s]$ for every $s\in I$ with $t_0<s<t_0+\varepsilon$.

Let $\beta=\displaystyle\sup\mathcal S$. For each $s\in\mathcal S$, denote by $\mathbf{u}_s$ the unique parallel section on $[t_0,s]$ with $\mathbf{u}_s(t_0)=v$. If $s_1,s_2\in\mathcal S$ and $s_1<s_2$, then $\mathbf{u}_{s_2}|_{[t_0,s_1]}$ is a parallel section on $[t_0,s_1]$ with initial value $v$. By the uniqueness defining $s_1\in\mathcal S$, we have $\mathbf{u}_{s_2}|_{[t_0,s_1]}=\mathbf{u}_{s_1}$. Thus the solutions $\mathbf{u}_s$ are mutually compatible.

Using this compatibility, define a section $\mathbf{u}$ on $[t_0,\beta)$ as follows. If $t\in[t_0,\beta)$, choose $s\in\mathcal S$ such that $t\leq s<\beta$ and set $\mathbf{u}(t)=\mathbf{u}_s(t)$. The definition is independent of the choice of $s$ precisely because of the preceding compatibility. Moreover, $\mathbf{u}$ is smooth and parallel on $[t_0,\beta)$ because it locally agrees with one of the smooth parallel sections $\mathbf{u}_s$.

Let us show that $\beta=b$. Suppose, for a contradiction, that $\beta<b$. Then $\beta\in I$. Choose a local trivialization of $\mathbf{E}$ over a neighborhood $U$ of $\gamma(\beta)$. Since $\gamma^{-1}(U)$ is open in $I$ and contains $\beta$, there exists $\delta>0$ such that $\gamma((\beta-\delta,\beta+\delta)\cap I)\subseteq U$. Choose $\tau\in\mathcal S$ with $\beta-\frac{\delta}{2}<\tau<\beta$. By the local case, there is a unique parallel section $\widetilde{\mathbf{u}}$ along $\gamma$ on $(\beta-\delta,\beta+\delta)\cap I$ such that $\widetilde{\mathbf{u}}(\tau)=\mathbf{u}(\tau)$.

We claim that $\widetilde{\mathbf{u}}=\mathbf{u}$ on $[\tau,\beta)$. Let $t\in[\tau,\beta)$. Choose $s\in\mathcal S$ such that $t\leq s<\beta$. On the interval $[\tau,s]$, both $\widetilde{\mathbf{u}}$ and $\mathbf{u}_s$ are parallel sections with the same value at $\tau$. By local uniqueness within the trivialization over $U$, we have $\widetilde{\mathbf{u}}=\mathbf{u}_s$ on $[\tau,s]$. Since $\mathbf{u}=\mathbf{u}_s$ on $[\tau,s]$, we conclude that $\widetilde{\mathbf{u}}(t)=\mathbf{u}(t)$. Since $t$ was arbitrary, $\widetilde{\mathbf{u}}=\mathbf{u}$ on $[\tau,\beta)$.

Now choose $\varepsilon>0$ such that $\beta+\varepsilon\in I$ and $\varepsilon<\delta$. Define a section $\widehat{\mathbf{u}}$ on $[t_0,\beta+\varepsilon]$ by
\[
\widehat{\mathbf{u}}(t)
=
\begin{cases}
\mathbf{u}(t), & t\in[t_0,\tau],\\
\widetilde{\mathbf{u}}(t), & t\in[\tau,\beta+\varepsilon].
\end{cases}
\]
This definition is smooth because $\mathbf{u}$ and $\widetilde{\mathbf{u}}$ agree on an interval to the right of $\tau$. Moreover, $\widehat{\mathbf{u}}$ is parallel because it locally agrees with parallel sections.

We also prove uniqueness of $\widehat{\mathbf{u}}$. Let $\mathbf{w}$ be another parallel section on $[t_0,\beta+\varepsilon]$ with $\mathbf{w}(t_0)=v$.
Then $\mathbf{w}|_{[t_0,\tau]}$ is a parallel section with the same initial value as $\mathbf{u}_\tau$. By uniqueness on $[t_0,\tau]$, we have
$\mathbf{w}=\mathbf{u}$ on $[t_0,\tau]$. In particular,
$\mathbf{w}(\tau)=\mathbf{u}(\tau)=\widetilde{\mathbf{u}}(\tau)$. By local uniqueness in the trivialization over $U$, we obtain
$\mathbf{w}=\widetilde{\mathbf{u}}$ on $[\tau,\beta+\varepsilon]$. Therefore,
$\mathbf{w}=\widehat{\mathbf{u}}$ on all of $[t_0,\beta+\varepsilon]$.

This shows that $\beta+\varepsilon\in\mathcal S$, contradicting the definition of $\beta=\displaystyle\sup\mathcal S$. Therefore $\beta=b$, and there is a unique parallel section to the right of $t_0$.

The same argument to the left, using intervals of the form $[s,t_0]$ with $s<t_0$, shows that there is a unique parallel section to the left of $t_0$. Gluing the two sections gives a unique parallel section on all of $I$ with initial value $v$.
\end{proof}

\begin{definition}[Parallel transport]\label{def:variedades-riemannianas-transporte-paralelo}\index{parallel transport}
Let $M$ be a smooth manifold with or without boundary, let
$\pi\colon \mathbf{E}\longrightarrow M$ be a smooth vector bundle with connection
$\nabla^{\mathbf{E}}$, and let $\gamma\colon I\longrightarrow M$ be a smooth curve. Given $t_0,t_1\in I$, define parallel transport from $t_0$ to
$t_1$ along $\gamma$ to be the map
\[
P^{\mathbf{E}}_{\gamma,t_0\to t_1}\colon \mathbf{E}_{\gamma(t_0)}\longrightarrow \mathbf{E}_{\gamma(t_1)}
\]
given by $P^{\mathbf{E}}_{\gamma,t_0\to t_1}(v)=\mathbf{u}(t_1)$, where
$\mathbf{u}$ is the unique parallel section along $\gamma$ such that
$\mathbf{u}(t_0)=v$.

When the bundle is clear from context, we omit the superscript and write $P_{\gamma,t_0\to t_1}$.
\end{definition}

\begin{proposition}\label{prop: transporte paralelo isomorfismo}
Let $M$ be a smooth manifold with or without boundary, let
$\pi\colon \mathbf{E}\longrightarrow M$ be a smooth vector bundle with connection
$\nabla^{\mathbf{E}}$, and let $\gamma\colon I\longrightarrow M$ be a smooth curve. For any $t_0,t_1\in I$, parallel transport
$P^{\mathbf{E}}_{\gamma,t_0\to t_1}\colon \mathbf{E}_{\gamma(t_0)}
\longrightarrow \mathbf{E}_{\gamma(t_1)}$ is a linear isomorphism. Moreover,
$P^{\mathbf{E}}_{\gamma,t_1\to t_0}=(P^{\mathbf{E}}_{\gamma,t_0\to t_1})^{-1}$.
More generally, if $t_0,t_1,t_2\in I$, then
\[
P^{\mathbf{E}}_{\gamma,t_1\to t_2}\circ P^{\mathbf{E}}_{\gamma,t_0\to t_1}
=
P^{\mathbf{E}}_{\gamma,t_0\to t_2}.
\]
\end{proposition}

\begin{proof}
Linearity follows from linearity of the parallelism equation. If
$\mathbf{u}$ and $\mathbf{w}$ are the parallel sections with initial conditions $v$ and $z$ at $t_0$, then $a\mathbf{u}+b\mathbf{w}$ is parallel and satisfies $(a\mathbf{u}+b\mathbf{w})(t_0)=av+bz$. By uniqueness,
$a\mathbf{u}+b\mathbf{w}$ is the parallel section corresponding to the initial data $av+bz$. Evaluating at $t_1$ gives
\[
P^{\mathbf{E}}_{\gamma,t_0\to t_1}(av+bz)
=
aP^{\mathbf{E}}_{\gamma,t_0\to t_1}(v)+bP^{\mathbf{E}}_{\gamma,t_0\to t_1}(z).
\]

To prove the inverse formula, let $v\in \mathbf{E}_{\gamma(t_0)}$ and let $\mathbf{u}$ be the parallel section with $\mathbf{u}(t_0)=v$. Then
$\mathbf{u}(t_1)=P^{\mathbf{E}}_{\gamma,t_0\to t_1}(v)$. If we transport
$\mathbf{u}(t_1)$ from $t_1$ to $t_0$, the resulting parallel section must again be $\mathbf{u}$, by uniqueness. Therefore,
\[
P^{\mathbf{E}}_{\gamma,t_1\to t_0}(P^{\mathbf{E}}_{\gamma,t_0\to t_1}(v))=v.
\]
The same argument, interchanging $t_0$ and $t_1$, shows that composition in the other order is also the identity. Thus $P^{\mathbf{E}}_{\gamma,t_1\to t_0}=(P^{\mathbf{E}}_{\gamma,t_0\to t_1})^{-1}$.

Finally, let $\mathbf{u}$ be the parallel section with $\mathbf{u}(t_0)=v$.
Then $P^{\mathbf{E}}_{\gamma,t_0\to t_1}(v)=\mathbf{u}(t_1)$. Transporting this vector from $t_1$ to $t_2$ again gives the value $\mathbf{u}(t_2)$, by uniqueness. Therefore,
\[
P^{\mathbf{E}}_{\gamma,t_1\to t_2}(P^{\mathbf{E}}_{\gamma,t_0\to t_1}(v))
=
\mathbf{u}(t_2)
=
P^{\mathbf{E}}_{\gamma,t_0\to t_2}(v).
\]
Since $v$ was arbitrary, the composition identity follows.
\end{proof}

\begin{remark}[Parallel frames along a curve]\label{obs:variedades-riemannianas-marcos-paralelos-a-lo-largo-de-una-curva}
Let $(b_1,\dots,b_r)$ be a basis of $\mathbf{E}_{\gamma(t_0)}$. For each
$\alpha\in\{1,\dots,r\}$, let $\mathbf{E}_\alpha$ be the parallel section along $\gamma$ such that $\mathbf{E}_\alpha(t_0)=b_\alpha$. Since parallel transport is a linear isomorphism, for each $t\in I$ the family
$(\mathbf{E}_1(t),\dots,\mathbf{E}_r(t))$ is a basis of
$\mathbf{E}_{\gamma(t)}$. This family of parallel sections is called a parallel frame along $\gamma$.

If $\mathbf{u}\in\Gamma^1(\gamma^*\mathbf{E})$, then there is a unique family of functions $u^1,\dots,u^r$ of class $C^1$ such that
$\mathbf{u}(t)=u^\alpha(t)\mathbf{E}_\alpha(t)$. Since
$\frac{D\mathbf{E}_\alpha}{dt}=0$, we obtain
\[
\displaystyle\frac{D\mathbf{u}}{dt}(t)=\displaystyle\frac{du^\alpha}{dt}(t)\mathbf{E}_\alpha(t).
\]
Consequently, a section $\mathbf{u}$ is parallel along $\gamma$ if and only if its components in a parallel frame are constant.
\end{remark}

Parallel transport along a single curve produces a frame along that curve. In a normal neighborhood, it can be performed simultaneously along all radial geodesics; this construction will provide the natural trivialization for spaces of sections on manifolds of bounded geometry.

\begin{definition}[Synchronous frame]
\label{def:marco-sincrono}
\index{synchronous frame}
Let $(M,\mathbf{g})$ be a Riemannian manifold without boundary, let $p\in M$, and let $U$ be a normal neighborhood of $p$. For $q\in U$, let
\[
\gamma_q(t):=\exp_p\bigl(t\exp_p^{-1}(q)\bigr),\qquad 0\leq t\leq1.
\]
If $\pi_{\mathbf{E}}\colon \mathbf{E}\longrightarrow M$ is a vector bundle with connection $\nabla^{\mathbf{E}}$ and
$(b_1,\dots,b_r)$ is a basis of $\mathbf{E}_p$, the frame
\[
\mathbf{e}_a(q):=P^{\mathbf{E}}_{\gamma_q,0\to1}(b_a),\qquad a\in\{1,\dots,r\},
\]
is called the \textit{synchronous frame centered at $p$}.
\end{definition}

\begin{proposition}[Properties of synchronous frames]
\label{prop:propiedades-marco-sincrono}
Let $(M,\mathbf{g})$ be a Riemannian manifold without boundary, let $p\in M$, let
$U$ be a normal neighborhood of $p$, and let
$\pi_{\mathbf{E}}\colon\mathbf{E}\longrightarrow M$ be a vector bundle with connection $\nabla^{\mathbf{E}}$. The synchronous frame centered at $p$ is smooth, and its restriction to every radial geodesic is parallel. If the connection is compatible with a bundle metric $\mathbf{h}_{\mathbf{E}}$ and the initial basis is orthonormal, the frame is orthonormal throughout $U$. In normal coordinates centered at $p$, if $\nabla^{\mathbf{E}}=d+\mathbf{A}$ and $\mathbf{A}=A_i\,\mathbf{d}x^i$, then
\[
\displaystyle\sum_{i=1}^n x^iA_i(x)=0,\qquad A_i(0)=0\quad(i\in\{1,\dots,n\}).
\]
\end{proposition}

\begin{proof}
Fix $q_0\in U$. The compact radial segment $\gamma_{q_0}([0,1])$
admits a finite subdivision $0=t_0<\cdots<t_N=1$ such that each
$\gamma_{q_0}([t_{j-1},t_j])$ lies in a trivializing coordinate domain.
By smooth dependence of $\gamma_q(t)=\exp_p(t\exp_p^{-1}q)$ on $(t,q)$
and compactness of each subinterval, the same domains and subdivision
work for all $q$ near $q_0$. On each subinterval parallel transport is
obtained from the linear system
\[
\frac{dv}{dt}(t,q)
=-\Gamma_i\bigl(\gamma_q(t)\bigr)\dot\gamma_q^i(t)v(t,q).
\]
Its solution depends smoothly on $q$ and on the initial vector. At the
subdivision points the change of frame is smooth as well. Composing these
finitely many transport maps proves smoothness of the synchronous frame
near $q_0$, and $q_0$ was arbitrary. Uniqueness of parallel transport
proves radial parallelism. Metric compatibility gives $\frac{d}{dt}\mathbf{h}_{\mathbf{E}}(\mathbf{e}_a,\mathbf{e}_b)=0$. In normal coordinates, the ray ending at $x$ is $t\mapsto tx$; parallelism implies
\[
0=\nabla^{\mathbf{E}}_{\displaystyle\sum_{i=1}^n x^i\boldsymbol{\partial}_i}\mathbf{e}_a
=\left(\displaystyle\sum_{i=1}^n x^iA_i(x)\right)\mathbf{e}_a.
\]
The radial identity follows because the $\mathbf{e}_a$ form a basis. Substituting
$x=tv$, dividing by $t>0$, and letting $t\to0$ gives $A_i(0)=0$.
\end{proof}

\begin{proposition}[Local equation for parallel transport]\label{prop: ecuacion local transporte paralelo operador}\index{local equation for parallel transport}
Let $M$ be a smooth manifold with or without boundary, let
$\pi\colon \mathbf{E}\longrightarrow M$ be a smooth vector bundle with connection
$\nabla^{\mathbf{E}}$, let $\gamma\colon I\longrightarrow M$ be a smooth curve, and suppose that $\gamma(J)\subseteq U$, where
$(\mathbf{e}_1,\dots,\mathbf{e}_r)$ is a local frame of $\mathbf{E}$ over
$U$. Fix $t_0\in J$, and represent parallel transport
$R(t):=P^{\mathbf{E}}_{\gamma,t_0\to t}\colon
\mathbf{E}_{\gamma(t_0)}\longrightarrow \mathbf{E}_{\gamma(t)}$ by a matrix $r(t)$ in the frame $\mathbf{e}_1,\dots,\mathbf{e}_r$. That is, if
$v=v^\beta \mathbf{e}_\beta(\gamma(t_0))$, then
\[
R(t)v=r^\alpha_{\ \beta}(t)v^\beta \mathbf{e}_\alpha(\gamma(t)).
\]
Then $r(t)$ satisfies
\[
\displaystyle\frac{dr}{dt}(t)+\dot\gamma^i(t)\Gamma_i(\gamma(t))r(t)=0,
\qquad
r(t_0)=I,
\]
where $\Gamma_i=(\Gamma^\alpha_{i\beta})_{\alpha,\beta}$. In particular, if $\dot\gamma(t_0)=\boldsymbol{\partial}_i|_{\gamma(t_0)}$, then $\displaystyle r'(t_0)=-\Gamma_i(\gamma(t_0))$.

Moreover, if $Q(t):=R(t)^{-1}=P^{\mathbf{E}}_{\gamma,t\to t_0}$ and $q(t)=r(t)^{-1}$ is its matrix, then $q(t)$ depends smoothly on $t$ and satisfies
\[
\displaystyle\frac{dq}{dt}(t)=q(t)\dot\gamma^i(t)\Gamma_i(\gamma(t)).
\]
\end{proposition}

\begin{proof}
Let $v=v^\beta \mathbf{e}_\beta(\gamma(t_0))$. By definition,
$\mathbf{u}(t):=R(t)v$ is the parallel section along $\gamma$ with
$\mathbf{u}(t_0)=v$. In the local frame, write
$\mathbf{u}(t)=u^\alpha(t)\mathbf{e}_\alpha(\gamma(t))$, where
$u^\alpha(t)=r^\alpha_{\ \beta}(t)v^\beta$.

Since $\mathbf{u}$ is parallel, the local equation
\eqref{eq: transporte paralelo haz} gives
\[
\displaystyle\frac{du^\alpha}{dt}(t)
=
-\dot\gamma^i(t)\Gamma^\alpha_{i\mu}(\gamma(t))u^\mu(t).
\]
Substituting $u^\alpha(t)=r^\alpha_{\ \beta}(t)v^\beta$, we obtain
\[
\displaystyle\frac{dr^\alpha_{\ \beta}}{dt}(t)v^\beta
=
-\dot\gamma^i(t)\Gamma^\alpha_{i\mu}(\gamma(t))r^\mu_{\ \beta}(t)v^\beta.
\]
Since $v$ is arbitrary, it follows that
\[
\displaystyle\frac{dr^\alpha_{\ \beta}}{dt}(t)
+
\dot\gamma^i(t)\Gamma^\alpha_{i\mu}(\gamma(t))r^\mu_{\ \beta}(t)
=
0.
\]
In matrix form, this is $\displaystyle \frac{dr}{dt}(t)+\dot\gamma^i(t)\Gamma_i(\gamma(t))r(t)=0$. Moreover, since $R(t_0)$ is the identity on $\mathbf{E}_{\gamma(t_0)}$, we have $r(t_0)=I$. If $\dot\gamma(t_0)=\boldsymbol{\partial}_i|_{\gamma(t_0)}$, then $\dot\gamma^j(t_0)=\delta_i^j$, and evaluating the equation at $t_0$ gives $r'(t_0)=-\Gamma_i(\gamma(t_0))$.

Finally, $q(t)=r(t)^{-1}$ is smooth because $r(t)$ is smooth and takes values in $GL(r,\mathbb K)$. Differentiating the identity $q(t)r(t)=I$ gives $q'(t)r(t)+q(t)r'(t)=0$. Therefore $q'(t)=-q(t)r'(t)r(t)^{-1}$. Using $r'(t)=-\dot\gamma^i(t)\Gamma_i(\gamma(t))r(t)$ yields $\displaystyle q'(t)=q(t)\dot\gamma^i(t)\Gamma_i(\gamma(t))$.
\end{proof}

\begin{theorem}[Covariant derivative through parallel transport]\label{teo: conexion sobre curvas y transporte paralelo}\index{covariant derivative through parallel transport}
Let $M$ be a smooth manifold with or without boundary, let
$\pi\colon \mathbf{E}\longrightarrow M$ be a smooth vector bundle with connection
$\nabla^{\mathbf{E}}$, let $\gamma\colon I\longrightarrow M$ be a smooth curve, and let $\mathbf{u}\in\Gamma^1(\gamma^*\mathbf{E})$. Then, for every
$t_0\in I$,
\[
\displaystyle\frac{D\mathbf{u}}{dt}(t_0)
=
\lim_{t\to t_0}
\displaystyle\frac{P^{\mathbf{E}}_{\gamma,t\to t_0}\mathbf{u}(t)-\mathbf{u}(t_0)}{t-t_0}.
\]
The limit is taken in the fixed vector space $\mathbf{E}_{\gamma(t_0)}$. If $t_0$ is an endpoint of the interval of definition, the limit is understood to be one-sided.
\end{theorem}

\begin{proof}
Choose a parallel frame $(\mathbf{E}_1,\dots,\mathbf{E}_r)$ along $\gamma$. Write $\mathbf{u}(t)=u^\alpha(t)\mathbf{E}_\alpha(t)$. By the preceding remark,
\[
\displaystyle\frac{D\mathbf{u}}{dt}(t_0)
=
\displaystyle\frac{du^\alpha}{dt}(t_0)\mathbf{E}_\alpha(t_0).
\]
On the other hand, since $\mathbf{E}_\alpha$ is parallel, parallel transport of $\mathbf{E}_\alpha(t)$ from $t$ to $t_0$ is $\mathbf{E}_\alpha(t_0)$. By linearity of parallel transport,
\[
P^{\mathbf{E}}_{\gamma,t\to t_0}\mathbf{u}(t)
=
u^\alpha(t)\mathbf{E}_\alpha(t_0).
\]
Then
\[
\displaystyle\frac{P^{\mathbf{E}}_{\gamma,t\to t_0}\mathbf{u}(t)-\mathbf{u}(t_0)}{t-t_0}
=
\displaystyle\frac{u^\alpha(t)-u^\alpha(t_0)}{t-t_0}\mathbf{E}_\alpha(t_0).
\]
Taking the limit as $t\to t_0$ gives
\[
\lim_{t\to t_0}
\displaystyle\frac{P^{\mathbf{E}}_{\gamma,t\to t_0}\mathbf{u}(t)-\mathbf{u}(t_0)}{t-t_0}
=
\displaystyle\frac{du^\alpha}{dt}(t_0)\mathbf{E}_\alpha(t_0)
=
\displaystyle\frac{D\mathbf{u}}{dt}(t_0).
\]
\end{proof}

\begin{corollary}\label{cor: conexion y transporte paralelo}
Let $M$ be a smooth manifold with or without boundary, and let
$\pi\colon \mathbf{E}\longrightarrow M$ be a smooth vector bundle with connection
$\nabla^{\mathbf{E}}$. Let $\mathbf{s}\in\Gamma^1(\mathbf{E})$, let $p\in M$, and let $\mathbf{X}\in\mathfrak X(M)$. If $\gamma$ is a smooth curve defined near $0$ such that $\gamma(0)=p$ and $\dot\gamma(0)=\mathbf{X}_p$, then
\[
\nabla^{\mathbf{E}}_{\mathbf{X}}\mathbf{s}(p)
=
\lim_{h\to 0}
\displaystyle\frac{P^{\mathbf{E}}_{\gamma,h\to 0}\mathbf{s}(\gamma(h))-\mathbf{s}(p)}{h}.
\]
\end{corollary}

\begin{proof}
Define a section $\mathbf{u}$ of $\gamma^*\mathbf{E}$ by
$\mathbf{u}(h)=\mathbf{s}(\gamma(h))$. Since $\mathbf{u}$ is the restriction of the section $\mathbf{s}$ to the curve $\gamma$,
Theorem~\ref{teo: derivada covariante a lo largo de una curva} implies that
\[
\displaystyle\frac{D\mathbf{u}}{dt}(0)
=
\nabla^{\mathbf{E}}_{\dot\gamma(0)}\mathbf{s}
=
\nabla^{\mathbf{E}}_{\mathbf{X}_p}\mathbf{s}
=
\nabla^{\mathbf{E}}_{\mathbf{X}}\mathbf{s}(p).
\]
Using Theorem~\ref{teo: conexion sobre curvas y transporte paralelo} at $t_0=0$, we obtain
\[
\displaystyle\frac{D\mathbf{u}}{dt}(0)
=
\lim_{h\to 0}
\displaystyle\frac{P^{\mathbf{E}}_{\gamma,h\to 0}\mathbf{u}(h)-\mathbf{u}(0)}{h}.
\]
Since $\mathbf{u}(h)=\mathbf{s}(\gamma(h))$ and $\mathbf{u}(0)=\mathbf{s}(p)$, the formula follows.
\end{proof}

\begin{lemma}[Fundamental theorem of covariant calculus]\label{lema: TFC covariante}\index{fundamental theorem of covariant calculus}
Let $M$ be a smooth manifold with or without boundary, and let
$\pi_{\mathbf{F}}\colon \mathbf{F}\longrightarrow M$ be a smooth vector bundle equipped with a connection $\nabla^{\mathbf{F}}$. Let
$\gamma\colon [a,b]\longrightarrow M$ be a smooth curve, and let
$\mathbf{s}\in\Gamma^1(\gamma^*\mathbf{F})$. For $t,t_0\in[a,b]$, denote by $P_{\gamma,t\to t_0}^{\mathbf{F}}\colon
\mathbf{F}_{\gamma(t)}\longrightarrow \mathbf{F}_{\gamma(t_0)}$ parallel transport along $\gamma$ from $\gamma(t)$ to
$\gamma(t_0)$. Then
\[
P_{\gamma,t\to t_0}^{\mathbf{F}}(\mathbf{s}(t))-\mathbf{s}(t_0)
=
\int_{t_0}^{t}
P_{\gamma,\tau\to t_0}^{\mathbf{F}}
\left(
\displaystyle\frac{D\mathbf{s}}{d\tau}(\tau)
\right)
\,d\tau.
\]
The integral is taken in the finite-dimensional vector space $\mathbf{F}_{\gamma(t_0)}$.
\end{lemma}

\begin{proof}
Fix $t_0\in[a,b]$. Choose a basis $v_1,\dots,v_r$ of the fiber $\mathbf{F}_{\gamma(t_0)}$. For each $\mu\in\{1,\dots,r\}$, let $\mathbf{e}_\mu$ be the parallel section of $\mathbf{F}$ along $\gamma$ satisfying $\mathbf{e}_\mu(t_0)=v_\mu$. Then $\mathbf{e}_1(t),\dots,\mathbf{e}_r(t)$ is a basis of $\mathbf{F}_{\gamma(t)}$ for every $t\in[a,b]$.

Since $\mathbf{s}$ is a section of $\mathbf{F}$ along $\gamma$, there exist functions $a^1,\dots,a^r\in C^1([a,b])$ such that $\displaystyle \mathbf{s}(t)=\displaystyle\sum_{\mu=1}^{r}a^\mu(t)\mathbf{e}_\mu(t)$. Each $\mathbf{e}_\mu$ is parallel, so its covariant derivative along $\gamma$ vanishes, and therefore $\displaystyle \frac{D\mathbf{s}}{dt}(t)=\displaystyle\sum_{\mu=1}^{r}\frac{d a^\mu}{dt}(t)\mathbf{e}_\mu(t)$.

Moreover, by the definition of parallel transport, $P_{\gamma,t\to t_0}^{\mathbf{F}}(\mathbf{e}_\mu(t))=\mathbf{e}_\mu(t_0)=v_\mu$. Consequently, $\displaystyle P_{\gamma,t\to t_0}^{\mathbf{F}}(\mathbf{s}(t))=\displaystyle\sum_{\mu=1}^{r}a^\mu(t)v_\mu$ and $\displaystyle P_{\gamma,t\to t_0}^{\mathbf{F}}\left(\frac{D\mathbf{s}}{dt}(t)\right)=\displaystyle\sum_{\mu=1}^{r}\frac{d a^\mu}{dt}(t)v_\mu$.

Applying the fundamental theorem of calculus componentwise in the basis $(v_1,\dots,v_r)$ of $\mathbf{F}_{\gamma(t_0)}$, we obtain
\[
P_{\gamma,t\to t_0}^{\mathbf{F}}(\mathbf{s}(t))-\mathbf{s}(t_0)
=
\displaystyle\sum_{\mu=1}^{r}\left(a^\mu(t)-a^\mu(t_0)\right)v_\mu.
\]
Since $\displaystyle a^\mu(t)-a^\mu(t_0)=\int_{t_0}^{t}\frac{d a^\mu}{d\tau}(\tau)\,d\tau$, it follows that
\[
P_{\gamma,t\to t_0}^{\mathbf{F}}(\mathbf{s}(t))-\mathbf{s}(t_0)
=
\int_{t_0}^{t}P_{\gamma,\tau\to t_0}^{\mathbf{F}}\left(\displaystyle\frac{D\mathbf{s}}{d\tau}(\tau)\right)\,d\tau,
\]
which proves the identity.
\end{proof}

\begin{lemma}[Integral criterion for the covariant derivative]\label{lema: criterio integral derivada covariante}\index{integral criterion for the covariant derivative}
Let $M$ be a smooth manifold with or without boundary, and let
$\pi_{\mathbf{F}}\colon \mathbf{F}\longrightarrow M$ be a smooth vector bundle equipped with a connection $\nabla^{\mathbf{F}}$. Let
$\mathbf{s}\in\Gamma^0(\mathbf{F})$ and
$\mathbf{A}\in\Gamma^0(T^*M\otimes \mathbf{F})$.

Suppose that, for every smooth curve $\gamma\colon [a,b]\longrightarrow M$ and any $t,t_0\in[a,b]$, we have
\[
P_{\gamma,t\to t_0}^{\mathbf{F}}(\mathbf{s}(\gamma(t)))-\mathbf{s}(\gamma(t_0))
=
\int_{t_0}^{t}
P_{\gamma,\tau\to t_0}^{\mathbf{F}}
\left(
\mathbf{A}_{\gamma(\tau)}(\dot\gamma(\tau))
\right)
\,d\tau.
\]
Then $\mathbf{s}\in\Gamma^1(\mathbf{F})$ and $\nabla^{\mathbf{F}} \mathbf{s}=\mathbf{A}$.
\end{lemma}

\begin{proof}
Fix $p\in M$. Choose a chart $(U,\varphi)$ around $p$, with $\varphi=(x^1,\dots,x^d)$, and a smooth local frame $\mathbf{e}_1,\dots,\mathbf{e}_r$ of $\mathbf{F}$ over $U$.

Write $\mathbf{s}=\displaystyle\sum_{\mu=1}^{r}s^\mu \mathbf{e}_\mu$, where the functions $s^\mu\colon U\longrightarrow\mathbb K$ are continuous. Likewise, write
\[
\mathbf{A}
=
\displaystyle\sum_{i=1}^{d}\displaystyle\sum_{\mu=1}^{r}
A_i^\mu\,\mathbf{d}x^i\otimes \mathbf{e}_\mu,
\]
where the functions $A_i^\mu\colon U\longrightarrow\mathbb K$ are continuous.

Denote by $\Gamma_{i\nu}^{\mu}$ the coefficients of the connection $\nabla^{\mathbf{F}}$ in the chosen frame, that is,
\[
\nabla^{\mathbf{F}}_{\boldsymbol{\partial}_i}\mathbf{e}_\nu
=
\displaystyle\sum_{\mu=1}^{r}\Gamma_{i\nu}^{\mu}\mathbf{e}_\mu,
\qquad
\boldsymbol{\partial}_i:=\displaystyle\boldsymbol{\partial}_{i}.
\]

Fix $i\in\{1,\dots,d\}$. We work with the coordinate curve in the direction $i$. If $p$ is an interior point, choose $\delta>0$ so that
\[
\gamma_i\colon (-\delta,\delta)\longrightarrow U,
\qquad
\gamma_i(t):=\varphi^{-1}\bigl(\varphi(p)+te_i\bigr),
\]
is well defined. If $p\in\partial M$ and a coordinate direction is defined only on one side within the boundary chart, take the corresponding one-sided interval. The following argument applies unchanged, replacing ordinary limits by one-sided limits in that direction. In all cases, $\gamma_i(0)=p$ and $\dot\gamma_i(t)=\boldsymbol{\partial}_i|_{\gamma_i(t)}$ for the admissible values of $t$.

Let $R_i(t)\colon \mathbf{F}_p\longrightarrow \mathbf{F}_{\gamma_i(t)}$ be parallel transport from $p=\gamma_i(0)$ to $\gamma_i(t)$, and let $Q_i(t):=R_i(t)^{-1}=P_{\gamma_i,t\to0}^{\mathbf{F}}$. In the local frame $\mathbf{e}_1,\dots,\mathbf{e}_r$, represent $R_i(t)$ and $Q_i(t)$ by matrices denoted by $r_i(t)$ and $q_i(t)$, respectively. Thus $q_i(t)=r_i(t)^{-1}$ and $r_i(0)=q_i(0)=I$.

By Proposition~\ref{prop: ecuacion local transporte paralelo operador}, the local differential equation for parallel transport gives
\[
r_i'(t)+\Gamma_i(\gamma_i(t))\,r_i(t)=0,
\]
where $\Gamma_i:=\bigl(\Gamma_{i\nu}^{\mu}\bigr)_{\mu,\nu}$. In particular, $r_i'(0)=-\Gamma_i(p)$.

Use the hypothesis of the lemma with the curve $\gamma_i$ and with $t_0=0$. Since $\dot\gamma_i(t)=\boldsymbol{\partial}_i|_{\gamma_i(t)}$, we obtain
\[
P_{\gamma_i,t\to0}^{\mathbf{F}}(\mathbf{s}(\gamma_i(t)))-\mathbf{s}(p)
=
\int_0^t
P_{\gamma_i,\tau\to0}^{\mathbf{F}}
\left(
\mathbf{A}_{\gamma_i(\tau)}(\boldsymbol{\partial}_i)
\right)
\,d\tau.
\]

Express this identity in the local frame. Define
\[
\widehat{s}:=
\begin{pmatrix}
s^1\\
\vdots\\
s^r
\end{pmatrix},
\qquad
\widehat{A}_i:=
\begin{pmatrix}
A_i^1\\
\vdots\\
A_i^r
\end{pmatrix}.
\]
Then
\[
q_i(t)\,\widehat{s}(\gamma_i(t))-\widehat{s}(p)
=
\int_0^t
q_i(\tau)\,\widehat{A}_i(\gamma_i(\tau))
\,d\tau.
\]
Multiplying by $r_i(t)$ gives
\[
\widehat{s}(\gamma_i(t))
=
r_i(t)\widehat{s}(p)
+
r_i(t)
\int_0^t
q_i(\tau)\,\widehat{A}_i(\gamma_i(\tau))
\,d\tau.
\]
Subtract $\widehat{s}(p)$ and divide by $t\neq0$:
\[
\displaystyle\frac{\widehat{s}(\gamma_i(t))-\widehat{s}(p)}{t}
=
\displaystyle\frac{r_i(t)-I}{t}\,\widehat{s}(p)
+
r_i(t)\,
\displaystyle\frac{1}{t}
\int_0^t
q_i(\tau)\,\widehat{A}_i(\gamma_i(\tau))
\,d\tau.
\]

Let $t\to0$ through the admissible parameter values of the curve. Since $r_i$ is differentiable and $r_i'(0)=-\Gamma_i(p)$, we have
\[
\displaystyle\frac{r_i(t)-I}{t}\longrightarrow -\Gamma_i(p).
\]
Moreover, $r_i(t)\to I$, $q_i(t)\to I$, and $\widehat{A}_i(\gamma_i(t))\to\widehat{A}_i(p)$ by continuity. Therefore the function $\tau\mapsto q_i(\tau)\widehat{A}_i(\gamma_i(\tau))$ is continuous at $0$, and consequently
\[
\displaystyle\frac{1}{t}
\int_0^t
q_i(\tau)\,\widehat{A}_i(\gamma_i(\tau))
\,d\tau
\longrightarrow
\widehat{A}_i(p).
\]
We conclude that
\[
\displaystyle\frac{\widehat{s}(\gamma_i(t))-\widehat{s}(p)}{t}
\longrightarrow
-\Gamma_i(p)\widehat{s}(p)+\widehat{A}_i(p).
\]

The left-hand side is the difference quotient of the components of $\mathbf{s}$ in the coordinate direction $\boldsymbol{\partial}_i$. Therefore $\partial_i\widehat{s}(p)$ exists, in the usual sense if $p$ is interior and in the corresponding one-sided sense if $p$ lies in a boundary chart, and satisfies
\[
\partial_i\widehat{s}(p)=\widehat{A}_i(p)-\Gamma_i(p)\widehat{s}(p).
\]

Since $p\in U$ was arbitrary, we obtain the identity $\partial_i\widehat{s}=\widehat{A}_i-\Gamma_i\widehat{s}$ throughout $U$. The right-hand side is continuous because $\widehat{A}_i$, $\Gamma_i$, and $\widehat{s}$ are continuous. Thus the partial derivatives of the components $s^\mu$ exist and are continuous in the corresponding charts. Therefore $\mathbf{s}$ is of class $C^1$ on $U$.

Finally, the local formula for the covariant derivative gives
\[
\nabla^{\mathbf{F}}_{\boldsymbol{\partial}_i}\mathbf{s}
=
\displaystyle\sum_{\mu=1}^{r}
\left(
\partial_i s^\mu
+
\displaystyle\sum_{\nu=1}^{r}
\Gamma_{i\nu}^{\mu}s^\nu
\right)\mathbf{e}_\mu.
\]
The identity $\partial_i\widehat{s}=\widehat{A}_i-\Gamma_i\widehat{s}$ implies
\[
\nabla^{\mathbf{F}}_{\boldsymbol{\partial}_i}\mathbf{s}
=
\displaystyle\sum_{\mu=1}^{r}A_i^\mu \mathbf{e}_\mu
=
\mathbf{A}(\boldsymbol{\partial}_i).
\]
Since this holds for every $i\in\{1,\dots,d\}$, we conclude that $\nabla^{\mathbf{F}} \mathbf{s}=\mathbf{A}$ on $U$.

Since $p\in M$ was arbitrary, it follows that $\mathbf{s}\in\Gamma^1(\mathbf{F})$ and $\nabla^{\mathbf{F}} \mathbf{s}=\mathbf{A}$ on all of $M$.
\end{proof}

\begin{remark}\label{obs:variedades-riemannianas-caso-particular-transporte-paralelo-permite-comparar}
In the particular case $\mathbf{E}=TM$, parallel transport allows us to compare tangent vectors at different points of a curve. If $\nabla$ is a connection on $TM$, the velocity $\dot\gamma$ of a smooth curve $\gamma$ is a section of $\gamma^*TM$, so it makes sense to consider $\displaystyle \frac{D\dot\gamma}{dt}$. This leads to the notion of a geodesic: a curve whose velocity is parallel along itself, that is, a curve satisfying $\displaystyle \frac{D\dot\gamma}{dt}=0$.
\end{remark}

\subsection{Horizontal--vertical decomposition induced by a connection}
\label{subsec:descomposicion-horizontal-vertical-conexion}

Covariant differentiation along curves also describes the geometry of the total space of a vector bundle. The following construction is first carried out for an arbitrary bundle; the cotangent decomposition used later in symbolic calculus will be a specialization.

Let $M$ be a finite-dimensional smooth manifold with or without boundary, let
$\pi_{\mathbf{E}}\colon \mathbf{E}\longrightarrow M$ be a smooth vector bundle of finite rank over $\mathbb K\in\{\mathbb R,\mathbb C\}$, and let $\nabla^{\mathbf{E}}$ be a linear connection on $\mathbf{E}$. For $e\in \mathbf{E}_x$, define the \textbf{vertical space} at
$e$ by
\begin{equation}
V_e\mathbf{E}:=\ker(d\pi_{\mathbf{E}})_e\subseteq T_e\mathbf{E}.
\label{eq:espacio-vertical-haz-vectorial}
\end{equation}
If $w\in \mathbf{E}_x$, the curve $t\mapsto e+tw$ remains in $\mathbf{E}_x$ and determines the
\textbf{vertical lift}
\begin{equation}
w_e^{\mathrm V}:=
\left.\frac{d}{dt}\right|_{t=0}(e+tw)\in V_e\mathbf{E}.
\label{eq:levantamiento-vertical-haz-vectorial}
\end{equation}
Here, including in the complex case, the real tangent space is used; the complex structure of $\mathbf{E}_x$ canonically induces the corresponding complex structure on $V_e\mathbf{E}$.

\begin{proposition}[Vertical subbundle and vertical lift]
\label{prop:subhaz-vertical-levantamiento-vertical}
\index{vertical subbundle}
\index{vertical lift}
Let $M$ be a finite-dimensional smooth manifold with or without boundary, and let
$\pi_{\mathbf{E}}\colon\mathbf{E}\longrightarrow M$ be a smooth vector bundle of finite rank. The set
\[
V\mathbf{E}:=\coprod_{e\in \mathbf{E}}V_e\mathbf{E}
\]
is a smooth vector subbundle of $T\mathbf{E}\to \mathbf{E}$. For each
$e\in \mathbf{E}_x$, the map
\begin{equation}
\mathbf{E}_x\longrightarrow V_e\mathbf{E},
\qquad
w\longmapsto w_e^{\mathrm V},
\label{eq:isomorfismo-vertical-canonico}
\end{equation}
is a canonical linear isomorphism.
\end{proposition}

\begin{proof}
Let $U\subseteq M$ be a trivializing open set, and write a point of
$\mathbf{E}|_U$ as $(x,u)\in U\times\mathbb K^r$. The differential of the projection is
\[
(d\pi_{\mathbf{E}})_{(x,u)}(X,W)=X,
\qquad
(X,W)\in T_xM\times\mathbb K^r.
\]
Therefore,
\[
V_{(x,u)}\mathbf{E}=\{(0,W)\mid W\in\mathbb K^r\}.
\]
These local descriptions show that $V\mathbf{E}$ has constant rank and smooth trivializations. In the same coordinates, the curve
$t\mapsto(x,u+tw)$ gives
\[
w_{(x,u)}^{\mathrm V}=(0,w).
\]
Thus \eqref{eq:isomorfismo-vertical-canonico} is a linear isomorphism.
The definition in terms of addition and scalar multiplication in the fiber uses no trivialization; hence the isomorphism is canonical.
\end{proof}

The connection selects a complement to the vertical subbundle. To construct it without choosing coordinates, we associate to each tangent vector of the total space the covariant derivative of the fiber vector representing it.

\begin{definition}[Connection map and horizontal subbundle]
\label{def:mapeo-conexion-subhaz-horizontal}
\index{connection map}
\index{horizontal subbundle}
Let $M$ be a finite-dimensional smooth manifold with or without boundary, and let
$\mathbf{E}\longrightarrow M$ be a smooth vector bundle equipped with a connection $\nabla^{\mathbf{E}}$. Let $v\in T_e\mathbf{E}$, with
$e\in \mathbf{E}_x$. In local bundle coordinates associated with a frame
$(\mathbf{e}_1,\ldots,\mathbf{e}_r)$, write
\[
e=u^\beta\mathbf{e}_\beta(x),\qquad
v=X^i\boldsymbol{\partial}_{i}\Big|_e
 +W^\alpha\boldsymbol{\partial}_{u^\alpha}\Big|_e,
\qquad
\nabla^{\mathbf{E}}_{\boldsymbol{\partial}_i}\mathbf{e}_\beta
=\Gamma^\alpha_{i\beta}\mathbf{e}_\alpha.
\]
The \textbf{connection map} is defined by
\begin{equation}
K^{\nabla^{\mathbf{E}}}(v):=
\bigl(W^\alpha+\Gamma^\alpha_{i\beta}(x)X^iu^\beta\bigr)
\mathbf{e}_\alpha(x)\in \mathbf{E}_x.
\label{eq:definicion-mapeo-conexion}
\end{equation}
The transformation law for the connection coefficients makes this expression independent of the coordinates and frame. If
$\mathbf{s}(t)\in\mathbf{E}_{\gamma(t)}$ is a smooth curve with
$\mathbf{s}(0)=e$ and $\dot{\mathbf{s}}(0)=v$, then
$K^{\nabla^{\mathbf{E}}}(v)=\left.\frac{D\mathbf{s}}{dt}\right|_{t=0}$.
The \textbf{horizontal space} at $e$ is
\begin{equation}
H_e^{\nabla^{\mathbf{E}}}\mathbf{E}:=\ker K^{\nabla^{\mathbf{E}}}_e,
\label{eq:espacio-horizontal-haz-vectorial}
\end{equation}
where $K^{\nabla^{\mathbf{E}}}_e$ denotes the restriction to $T_e\mathbf{E}$.
\end{definition}

\begin{proposition}[Horizontal--vertical decomposition]
\label{prop:descomposicion-horizontal-vertical-conexion}
\index{horizontal lift}
Let $M$ be a finite-dimensional smooth manifold with or without boundary, and let
$\pi_{\mathbf{E}}\colon\mathbf{E}\longrightarrow M$ be a smooth vector bundle of finite rank equipped with a connection $\nabla^{\mathbf{E}}$. The local expressions
\eqref{eq:definicion-mapeo-conexion} are compatible on overlaps and define a smooth map
$K^{\nabla^{\mathbf{E}}}\colon T\mathbf{E}\to \mathbf{E}$ covering $\pi_{\mathbf{E}}$, and each
\[
K^{\nabla^{\mathbf{E}}}_e\colon T_e\mathbf{E}\longrightarrow \mathbf{E}_x
\]
is linear over $\mathbb R$; its restriction to the vertical subspace is linear over $\mathbb K$ under the isomorphism
\eqref{eq:isomorfismo-vertical-canonico}. Moreover,
\[
K^{\nabla^{\mathbf{E}}}_e(w_e^{\mathrm V})=w,
\qquad w\in \mathbf{E}_x,
\]
and we have the smooth decompositions
\begin{equation}
T_e\mathbf{E}=H_e^{\nabla^{\mathbf{E}}}\mathbf{E}\oplus V_e\mathbf{E},
\qquad
T\mathbf{E}=H^{\nabla^{\mathbf{E}}}\mathbf{E}\oplus V\mathbf{E}.
\label{eq:descomposicion-horizontal-vertical}
\end{equation}
The restriction
\begin{equation}
(d\pi_{\mathbf{E}})_e\big|_{H_e^{\nabla^{\mathbf{E}}}\mathbf{E}}\colon
H_e^{\nabla^{\mathbf{E}}}\mathbf{E}\longrightarrow T_xM
\label{eq:proyeccion-horizontal-isomorfismo}
\end{equation}
is an isomorphism. Its inverse sends $X\in T_xM$ to the
\textbf{horizontal lift} $X_e^{\mathrm H}$.
\end{proposition}

\begin{proof}
Choose coordinates $(x^1,\ldots,x^n)$ on an open set $U\subseteq M$ and a local frame $(\mathbf{e}_1,\ldots,\mathbf{e}_r)$ of $\mathbf{E}|_U$. Write
\[
\nabla^{\mathbf{E}}_{\boldsymbol{\partial}_i}\mathbf{e}_\beta
=\Gamma^\alpha_{i\beta}\mathbf{e}_\alpha.
\]
If $\mathbf{s}(t)=u^\alpha(t)\mathbf{e}_\alpha(\gamma(t))$, the local formula for covariant differentiation along a curve, obtained in
\eqref{eq: derivada covariante curva haz local}, gives
\begin{equation}
K^{\nabla^{\mathbf{E}}}(x,u;X,W)
=
\left(x,
W^\alpha+\Gamma^\alpha_{i\beta}(x)X^iu^\beta
\right).
\label{eq:mapeo-conexion-formula-local}
\end{equation}
The right-hand side depends only on $(x,u;X,W)=\dot{\mathbf{s}}(0)$ and not on the higher-order terms of the curve. This proves independence of the representing curve. The same formula proves smoothness and linearity of
$K^{\nabla^{\mathbf{E}}}_e$ in $(X,W)$.

For a vertical vector, $X=0$; therefore,
\[
K^{\nabla^{\mathbf{E}}}(x,u;0,W)=(x,W).
\]
By \eqref{eq:isomorfismo-vertical-canonico}, the restriction of
$K^{\nabla^{\mathbf{E}}}_e$ to $V_e\mathbf{E}$ is precisely the inverse of the vertical lift. In particular,
$H_e^{\nabla^{\mathbf{E}}}\mathbf{E}\cap V_e\mathbf{E}=\{0\}$. For
$v\in T_e\mathbf{E}$, set
\[
v^{\mathrm V}:=\bigl(K^{\nabla^{\mathbf{E}}}(v)\bigr)_e^{\mathrm V}.
\]
Then $v-v^{\mathrm V}\in H_e^{\nabla^{\mathbf{E}}}\mathbf{E}$ and
$v=(v-v^{\mathrm V})+v^{\mathrm V}$. This proves the first equality in
\eqref{eq:descomposicion-horizontal-vertical}. The local formula
\eqref{eq:mapeo-conexion-formula-local} also shows that the kernels
$H_e^{\nabla^{\mathbf{E}}}\mathbf{E}$ vary smoothly, giving the equality of subbundles.

If a horizontal vector belongs to the kernel of $(d\pi_{\mathbf{E}})_e$, it is also vertical and hence, by the direct sum, zero. The restriction
\eqref{eq:proyeccion-horizontal-isomorfismo} is therefore injective. For
$X=X^i\boldsymbol{\partial}_i|_x$, the formula
\begin{equation}
X_e^{\mathrm H}
=X^i\left(
\boldsymbol{\partial}_{i}
-\Gamma^\alpha_{i\beta}(x)u^\beta
\boldsymbol{\partial}_{u^\alpha}
\right)_{(x,u)}
\label{eq:levantamiento-horizontal-haz-local}
\end{equation}
defines a horizontal vector projecting to $X$. This proves surjectivity and completes the proof.
\end{proof}

A curve $\mathbf{s}(t)\in \mathbf{E}$ is horizontal if $\dot{\mathbf{s}}(t)$ belongs to
$H_{\mathbf{s}(t)}^{\nabla^{\mathbf{E}}}\mathbf{E}$ for every $t$. By the definition of the connection map, this is equivalent to $\frac{D\mathbf{s}}{dt}=0$; that is, $\mathbf{s}$ is horizontal if and only if the fiber vector it represents is parallel along
$\pi_{\mathbf{E}}\circ \mathbf{s}$. Consequently, $H^{\nabla^{\mathbf{E}}}\mathbf{E}$ is the Ehresmann connection induced by the linear connection $\nabla^{\mathbf{E}}$. The subbundle $V\mathbf{E}$ is canonical, whereas $H^{\nabla^{\mathbf{E}}}\mathbf{E}$ depends on the chosen connection. This construction should not be confused with the nonlinear connections introduced in Finsler geometry: here the horizontal distribution comes from a linear connection on a vector bundle.

We now specialize the construction to the cotangent bundle. Let $\nabla$ be a linear connection on $TM$, and equip $T^*M$ with the dual connection characterized in Proposition~\ref{conexion en tensores, duales y endomorfismos}. Let
$\pi\colon T^*M\to M$ be the projection. Then
\[
V(T^*M)=\ker d\pi,
\qquad
H^\nabla(T^*M):=H^{\nabla^{T^*M}}(T^*M).
\]

\begin{proposition}[Horizontal lift in the cotangent bundle]
\label{prop:levantamiento-horizontal-cotangente}
\index{horizontal lift!in the cotangent bundle}
Let $M$ be a finite-dimensional smooth manifold with or without boundary, let
$\nabla$ be a linear connection on $TM$, and let $(x^1,\ldots,x^n)$ be local coordinates. Write
\[
\nabla_{\boldsymbol{\partial}_i}\boldsymbol{\partial}_j
=\Gamma^k_{ij}\boldsymbol{\partial}_k,
\qquad
\boldsymbol{\xi}=\xi_j\,\mathbf{d}x^j.
\]
With these conventions, the horizontal lift of
$\boldsymbol{\partial}_i$ to
$T^*M$ is
\begin{equation}
\mathbf{H}_i
=\boldsymbol{\partial}_{i}
+\Gamma^k_{ij}(x)\xi_k
\frac{\partial}{\partial\xi_j}.
\label{eq:levantamiento-horizontal-cotangente-local}
\end{equation}
\end{proposition}

\begin{proof}
Let $x(t)$ be a curve, and let
$\boldsymbol{\xi}(t)=\xi_j(t)\mathbf{d}x^j|_{x(t)}$ be a section of the cotangent bundle along it.
The dual connection satisfies
\[
(\nabla_{\dot x}\boldsymbol{\xi})_j
=\dot\xi_j-\Gamma^k_{ij}(x(t))\dot x^i\xi_k.
\]
The curve $(x(t),\boldsymbol{\xi}(t))$ is horizontal if and only if the left-hand side vanishes. Therefore,
\[
\dot\xi_j
=\Gamma^k_{ij}(x(t))\dot x^i\xi_k.
\]
If $\dot x(0)=\boldsymbol{\partial}_i|_x$, the vertical component of the velocity is
$\Gamma^k_{ij}(x)\xi_k\boldsymbol{\partial}_{\xi_j}$, proving
\eqref{eq:levantamiento-horizontal-cotangente-local} and, in particular, fixing its sign. By contrast, for
$v=v^j\boldsymbol{\partial}_j\in T_xM$, the equation
$\nabla_{\dot x}v=0$ gives
$\dot v^k=-\Gamma^k_{ij}\dot x^iv^j$; according to
\eqref{eq:levantamiento-horizontal-haz-local}, the lift in $TM$ has the opposite sign.
\end{proof}

It remains to specify how this decomposition acts on sections of the pullback of a homomorphism bundle. Let $\mathbf{E},\mathbf{F}\to M$ be smooth vector bundles of finite rank with connections $\nabla^{\mathbf{E}}$ and $\nabla^{\mathbf{F}}$. The induced connection
$\nabla^{\operatorname{Hom}(\mathbf{E},\mathbf{F})}$ is characterized by
\eqref{eq:conexion-haz-homomorfismos}. The pullback connection on
$\pi^*\operatorname{Hom}(\mathbf{E},\mathbf{F})\to T^*M$ is the unique
connection that, for a local section $\mathbf{A}$ of
$\operatorname{Hom}(\mathbf{E},\mathbf{F})$ and
$\mathbf{Z}\in\mathfrak X(T^*M)$, satisfies
\begin{equation}
\begin{aligned}
&\left[
(\pi^*\nabla^{\operatorname{Hom}(\mathbf{E},\mathbf{F})})_{\mathbf{Z}}
(\pi^*\mathbf{A})
\right]_{(x,\xi)}\\
&\qquad=
\left(\nabla^{\operatorname{Hom}(\mathbf{E},\mathbf{F})}\mathbf{A}\right)_x
\bigl(d\pi_{(x,\xi)}(\mathbf{Z}_{(x,\xi)})\bigr).
\end{aligned}
\label{eq:conexion-pullback-hom}
\end{equation}
The equality is understood at each point $(x,\xi)$ lying over the domain
of $\mathbf{A}$, using the canonical identification
$(\pi^*\operatorname{Hom}(\mathbf{E},\mathbf{F}))_{(x,\xi)}
\cong\operatorname{Hom}(\mathbf{E}_x,\mathbf{F}_x)$.
The Leibniz rule determines the connection on all sections.
Existence and uniqueness are verified in the pullback frames associated
with the trivializations of
Proposition~\ref{prop: pullback haz vectorial suave}.
If $\omega$ is the local matrix of connection 1-forms of
$\nabla^{\operatorname{Hom}(\mathbf{E},\mathbf{F})}$, the matrix in the
pullback frame is $\pi^*\omega$.

For
\[
\mathbf{a}\in\Gamma\!\left(\pi^*\operatorname{Hom}(\mathbf{E},\mathbf{F})\right),
\quad X\in T_xM,
\quad \eta\in T_x^*M,
\quad \xi\in T_x^*M,
\]
define, at the point $(x,\xi)$,
\begin{equation}
\nabla_X^{\mathrm H}\mathbf{a}
:=(\pi^*\nabla^{\operatorname{Hom}(\mathbf{E},\mathbf{F})})_{X_{(x,\xi)}^{\mathrm H}}\mathbf{a},
\qquad
\nabla_\eta^{\mathrm V}\mathbf{a}
:=(\pi^*\nabla^{\operatorname{Hom}(\mathbf{E},\mathbf{F})})_{\eta_{(x,\xi)}^{\mathrm V}}\mathbf{a}.
\label{eq:definicion-derivadas-horizontal-vertical-cotangente}
\end{equation}
These are, respectively, the \textbf{horizontal covariant derivative} and the
\textbf{vertical covariant derivative} of $\mathbf{a}$.
\index{horizontal covariant derivative}
\index{vertical covariant derivative}
\index{pullback connection}

\begin{proposition}[Local formulas in the cotangent bundle]
\label{prop:formulas-locales-derivadas-horizontal-vertical}
Let $M$ be a smooth manifold with or without boundary. In the preceding coordinates, choose local frames of $\mathbf{E}$ and $\mathbf{F}$, and denote their connection matrices by $\omega_i^{\mathbf{E}}$ and $\omega_i^{\mathbf{F}}$. If
$a(x,\xi)$ is the local matrix of a section of
$\pi^*\operatorname{Hom}(\mathbf{E},\mathbf{F})$, then
\begin{equation}
\nabla_j^{\mathrm V}a
=\frac{\partial a}{\partial\xi_j},
\qquad
\nabla_i^{\mathrm H}a
=\frac{\partial a}{\partial x^i}
+\Gamma^k_{ij}\xi_k\frac{\partial a}{\partial\xi_j}
+\omega_i^{\mathbf{F}} a-a\omega_i^{\mathbf{E}}.
\label{eq:derivadas-geometricas-cotangente-locales}
\end{equation}
The two derivatives defined in
\eqref{eq:definicion-derivadas-horizontal-vertical-cotangente} are independent of the coordinates and frames used to calculate these expressions.
\end{proposition}

\begin{proof}
The vertical vector $\partial_{\xi_j}$ projects to zero. By
\eqref{eq:conexion-pullback-hom}, the pullback connection therefore contributes no connection term, leaving
$\nabla_j^{\mathrm V}a=\frac{\partial a}{\partial\xi_j}$. For the horizontal direction, \eqref{eq:levantamiento-horizontal-cotangente-local} differentiates the components of $a$ by
\[
\boldsymbol{\partial}_{i}
+\Gamma^k_{ij}\xi_k\frac{\partial}{\partial\xi_j}.
\]
The connection on $\operatorname{Hom}(\mathbf{E},\mathbf{F})$, given by
\eqref{eq:conexion-haz-homomorfismos}, adds exactly
$\omega_i^{\mathbf{F}}a-a\omega_i^{\mathbf{E}}$. This proves
\eqref{eq:derivadas-geometricas-cotangente-locales}. Finally,
$X^{\mathrm H}$ and $\eta^{\mathrm V}$ were defined as intrinsic lifts, and the pullback connection is also intrinsic. Thus a change of coordinates or frames merely transforms the matrices in the formula according to the laws for these objects and does not alter the resulting derivatives.
\end{proof}

\section{Geodesics}

Euclidean straight lines, parametrized with constant speed, satisfy $\gamma''(t)=0$ and realize the distance between their points. Distinguishing the image of a curve from its parametrization is essential: a nonaffine reparametrization of a straight line preserves length but generally has nonzero acceleration.

On a Riemannian manifold, it is natural to seek an analogue of these curves. The Euclidean notion of a straight line no longer makes sense directly, but the metric structure allows us to measure lengths and formulate the problem of finding curves that connect two points with minimal length. We will define \textit{geodesics} by their covariant acceleration. For the Levi--Civita connection, we will prove that geodesics minimize on sufficiently small intervals and that every constant-speed minimizing curve is a geodesic; an extended geodesic may cease to minimize.

\begin{figura}[H]
    \includegraphics[width=0.55\linewidth]{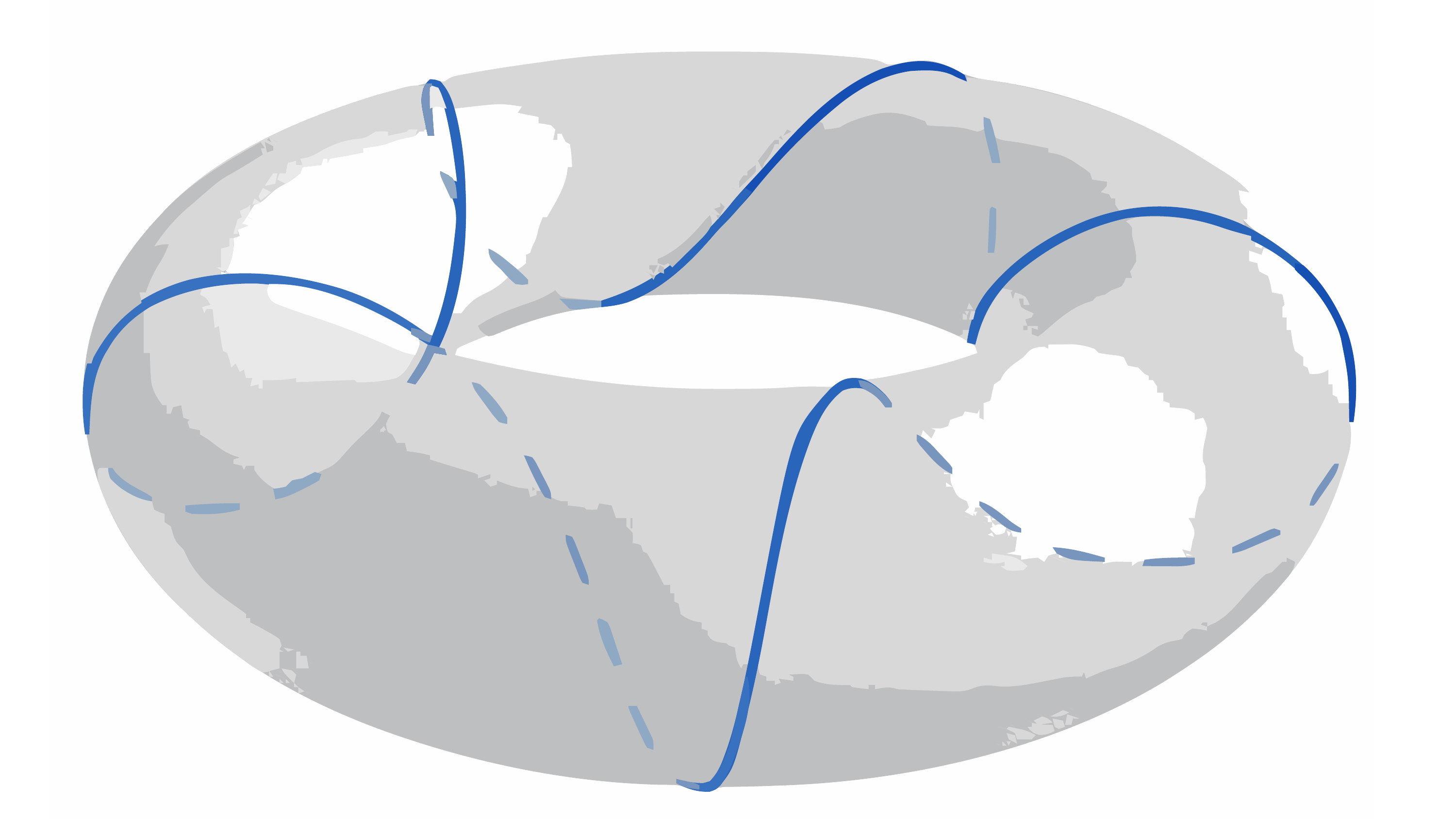}
    \caption{Geodesics on a torus. Their global behavior depends on the metric and initial direction; in general, being closed or dense does not exhaust the possibilities.}
    \label{fig:geodesicas-toro}
\end{figura}
However, to obtain an intrinsic characterization analogous to the condition $\gamma''(t)=0$ in $\mathbb{R}^n$, we cannot simply differentiate the curve as though all its velocities belonged to the same vector space. If $\gamma\colon I\longrightarrow M$ is a smooth curve, its velocity $\dot\gamma(t)$ belongs to the tangent space $T_{\gamma(t)}M$, which depends on the point $\gamma(t)$. Thus it makes no sense to compare $\dot\gamma(t)$ directly with $\dot\gamma(t+h)$ by ordinary subtraction.

The preceding section resolves precisely this difficulty: a connection on $TM$ induces a covariant derivative on sections of $\gamma^*TM$. Since the velocity $\dot\gamma$ is a section of $\gamma^*TM$, we may consider its covariant derivative $\frac{D\dot\gamma}{dt}$. This allows us to formulate the condition of vanishing covariant acceleration and thereby define the geodesics associated with a connection.

\begin{definition}\label{def:geodesica}\index{geodesic}
Let $M$ be a smooth manifold with or without boundary and let $\nabla$ be a connection on $TM$. We say that a smooth curve $\gamma\colon I\longrightarrow M$ is a \textbf{geodesic} if
\[
\displaystyle\frac{D\dot\gamma}{dt}\equiv 0.
\]

Equivalently, if $(U,\phi)$ is a smooth chart such that $\gamma(I)\subset U$ and
\[
x(t)=(\phi\circ\gamma)(t)=(x^{1}(t),\dots,x^{n}(t)),
\]
then $\gamma$ is a geodesic if and only if it satisfies the geodesic equations
\begin{equation}\label{eq: ecuacion de la geodesica}
\displaystyle\frac{d^{2}x^{k}}{dt^{2}}(t)
+
\displaystyle\frac{dx^{i}}{dt}(t)\displaystyle\frac{dx^{j}}{dt}(t)\Gamma_{ij}^{k}(x(t))
=
0.
\end{equation}
\end{definition}

We have defined geodesics, which in coordinates are described by a system of second-order ordinary differential equations. Before studying them in greater depth, we must ask whether such curves actually exist for every initial condition. Geodesics exist locally, as guaranteed by the fundamental theorem on flows. To see this geometrically, we first introduce the \textit{geodesic vector field}.

\begin{figura}[h]
    \includegraphics[width=0.7\linewidth]{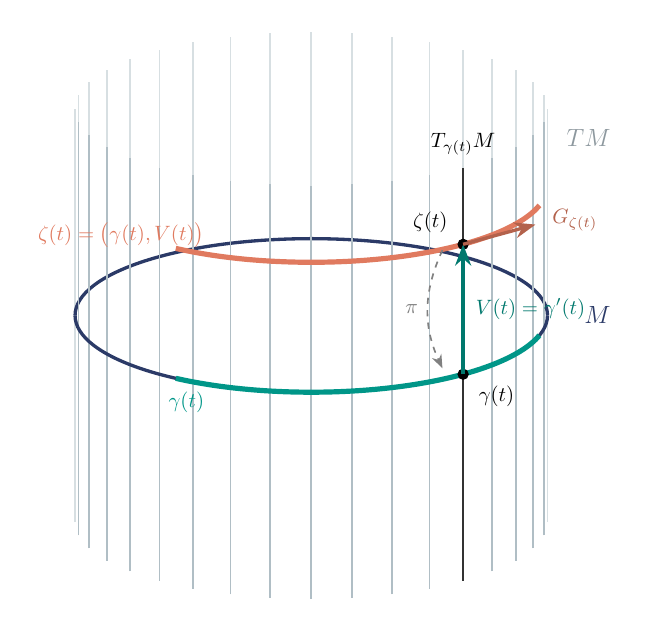}
    \caption{Schematic representation of the geodesic vector field on $TM$. The curve
    $\zeta(t)=(\gamma(t),\mathbf{V}(t))$ is an integral curve of $\mathbf{G}$
    if and only if $\mathbf{V}(t)=\dot{\gamma}(t)$ and
    $\frac{D\mathbf{V}}{dt}(t)=0$; in particular, when
    $\mathbf{V}(t)=\dot{\gamma}(t)$, the integral curves of $\mathbf{G}$
    correspond to the
    geodesics of $M$.}
    \label{fig:campo-geodesico}
\end{figura}

\begin{proposition}[Geodesic vector field]\label{prop: campo geodesico}\index{geodesic vector field}
Let $M$ be a smooth manifold without boundary and let $\nabla$ be a connection on $TM$. Then there exists a unique smooth vector field $\mathbf{G}\in \mathfrak{X}(TM)$ with the following property: for every smooth curve
\[
\zeta\colon I\longrightarrow TM,
\qquad
\zeta(t)=\bigl(\gamma(t),\mathbf{V}(t)\bigr),
\]
$\zeta$ is an integral curve of $\mathbf{G}$ if and only if
\[
\mathbf{V}(t)=\dot\gamma(t)
\quad\text{and}\quad
\displaystyle\frac{D\mathbf{V}}{dt}(t)=0
\]
for every $t\in I$.

Moreover, if $(U,\phi)$ is a smooth chart of $M$ and $(x^{1},\dots,x^{n},v^{1},\dots,v^{n})$ are the induced coordinates on $\pi^{-1}(U)\subset TM$, then $\mathbf{G}$ is given by
\[
\mathbf{G}
=
v^{k}\boldsymbol{\partial}_{x^k}
-
v^{i}v^{j}\Gamma_{ij}^{k}(x)\boldsymbol{\partial}_{v^k}.
\]

In particular, if $\zeta\colon I\longrightarrow TM$ is given by $\zeta(t)=\bigl(\gamma(t),\dot\gamma(t)\bigr)$, then $\zeta$ is an integral curve of $\mathbf{G}$ if and only if $\gamma$ is a geodesic.
\end{proposition}

\begin{proof}
Let $(U,\phi)$ be a smooth chart of $M$. On $\pi^{-1}(U)$, consider the vector field given, in the induced coordinates $(x^{1},\dots,x^{n},v^{1},\dots,v^{n})$, by
\[
\mathbf{G}_U
=
v^{k}\boldsymbol{\partial}_{x^k}
-
v^{i}v^{j}\Gamma_{ij}^{k}(x)\boldsymbol{\partial}_{v^k}.
\]

Let us first determine which curves are integral curves of $\mathbf{G}_U$. Let
$\zeta\colon I\longrightarrow \pi^{-1}(U)$ be a smooth curve and write
$\zeta(t)=\bigl(\gamma(t),\mathbf{V}(t)\bigr)$, where
$\gamma=\pi\circ \zeta$ and $\mathbf{V}(t)\in T_{\gamma(t)}M$ for every
$t\in I$. Let $\widetilde{\phi}\colon \pi^{-1}(U)\longrightarrow
\phi(U)\times \mathbb{R}^{n}$ be the diffeomorphism defining the induced
coordinates on $TM$. Write the coordinate representation of $\zeta$ as
$\widetilde{\phi}\circ \zeta=(\xi^{1},\dots,\xi^{n},\eta^{1},\dots,
\eta^{n})\colon I\longrightarrow \phi(U)\times \mathbb{R}^{n}$. Then
$\phi\circ\gamma=(\xi^{1},\dots,\xi^{n})$ and
$\mathbf{V}(t)=\eta^{k}(t)\boldsymbol{\partial}_k|_{\gamma(t)}$.

By the definition of an integral curve, $\zeta$ is an integral curve of $\mathbf{G}_U$ if and only if
\[
\displaystyle\frac{d\xi^{k}}{dt}(t)=\eta^{k}(t),
\qquad
\displaystyle\frac{d\eta^{k}}{dt}(t)
=
-\Gamma_{ij}^{k}(\xi(t))\eta^{i}(t)\eta^{j}(t)
\]
for every $t\in I$.

The first equation is equivalent to $\mathbf{V}(t)=\dot\gamma(t)$. Since
$\phi\circ \gamma=(\xi^{1},\dots,\xi^{n})$, we have
$\dot\gamma(t)=\frac{d\xi^{k}}{dt}(t)
\boldsymbol{\partial}_k|_{\gamma(t)}$, so
$\dot\gamma(t)=\mathbf{V}(t)$ if and only if
$\frac{d\xi^{k}}{dt}(t)=\eta^{k}(t)$ for every $k$.

On the other hand, by Proposition~\ref{prop: expresion local derivada covariante curvas},
\[
\displaystyle\frac{D\mathbf{V}}{dt}(t)
=
\left(
\displaystyle\frac{d\eta^{k}}{dt}(t)
+
\displaystyle\frac{d\xi^{i}}{dt}(t)\eta^{j}(t)\Gamma_{ij}^{k}(\xi(t))
\right)
\boldsymbol{\partial}_k\big|_{\gamma(t)}.
\]
If, in addition, $\mathbf{V}(t)=\dot\gamma(t)$, then
$\frac{d\xi^{i}}{dt}(t)=\eta^{i}(t)$, and consequently
$\frac{D\mathbf{V}}{dt}(t)=0$ is equivalent to
\[
\displaystyle\frac{d\eta^{k}}{dt}(t)
=
-\Gamma_{ij}^{k}(\xi(t))\eta^{i}(t)\eta^{j}(t).
\]

We have therefore proved that $\zeta$ is an integral curve of $\mathbf{G}_U$
if and only if $\mathbf{V}(t)=\dot\gamma(t)$ and
$\frac{D\mathbf{V}}{dt}(t)=0$ for every $t\in I$.

We now prove that these local vector fields agree on overlaps. Let $(U,\phi)$ and $(\widetilde{U},\widetilde{\phi})$ be two smooth charts with $U\cap \widetilde{U}\neq \varnothing$, and let $F=\widetilde{\phi}\circ \phi^{-1}$. If $(x^{1},\dots,x^{n},v^{1},\dots,v^{n})$ and $(\widetilde{x}^{1},\dots,\widetilde{x}^{n},\widetilde{v}^{1},\dots,\widetilde{v}^{n})$ denote the corresponding induced coordinates, then $\widetilde{x}^{a}=F^{a}(x)$ and $\widetilde{v}^{a}=\frac{\partial F^{a}}{\partial x^{i}}(x)v^{i}$.

Let $\zeta$ be an integral curve of $\mathbf{G}_U$ contained in $\pi^{-1}(U\cap \widetilde{U})$. If we write its coordinate representation in the system $(x,v)$ as $\widetilde{\phi}_U\circ \zeta=(\xi^{1},\dots,\xi^{n},\eta^{1},\dots,\eta^{n})$, then
\[
\displaystyle\frac{d\xi^{k}}{dt}(t)=\eta^{k}(t),
\qquad
\displaystyle\frac{d\eta^{k}}{dt}(t)
=
-\Gamma_{ij}^{k}(\xi(t))\eta^{i}(t)\eta^{j}(t).
\]

The components of the same curve in the system $(\widetilde{x},\widetilde{v})$ are $\widetilde{\xi}^{a}(t)=F^{a}(\xi(t))$ and $\widetilde{\eta}^{a}(t)=\frac{\partial F^{a}}{\partial x^{i}}(\xi(t))\eta^{i}(t)$. By the chain rule,
\[
\displaystyle\frac{d\widetilde{\xi}^{a}}{dt}(t)
=
\widetilde{\eta}^{a}(t).
\]

Moreover,
\[
\displaystyle\frac{d\widetilde{\eta}^{a}}{dt}(t)
=
\displaystyle\frac{\partial^{2}F^{a}}{\partial x^{i}\partial x^{j}}(\xi(t))\eta^{i}(t)\eta^{j}(t)
+
\displaystyle\frac{\partial F^{a}}{\partial x^{k}}(\xi(t))\displaystyle\frac{d\eta^{k}}{dt}(t).
\]
Substituting the equation for $\displaystyle\frac{d\eta^{k}}{dt}$ gives
\[
\displaystyle\frac{d\widetilde{\eta}^{a}}{dt}(t)
=
\left(
\displaystyle\frac{\partial^{2}F^{a}}{\partial x^{i}\partial x^{j}}
-
\displaystyle\frac{\partial F^{a}}{\partial x^{k}}\Gamma_{ij}^{k}
\right)(\xi(t))\eta^{i}(t)\eta^{j}(t).
\]
By the transformation law for the Christoffel symbols (Lemma~\ref{lema: transformacion metrica Christoffel}),
\[
\widetilde{\Gamma}_{bc}^{a}(F(x))
\displaystyle\frac{\partial F^{b}}{\partial x^{i}}(x)
\displaystyle\frac{\partial F^{c}}{\partial x^{j}}(x)
=
\displaystyle\frac{\partial F^{a}}{\partial x^{k}}(x)\Gamma_{ij}^{k}(x)
-
\displaystyle\frac{\partial^{2}F^{a}}{\partial x^{i}\partial x^{j}}(x).
\]
Applying this identity in $x=\xi(t)$, we conclude that
\[
\displaystyle\frac{d\widetilde{\eta}^{a}}{dt}(t)
=
-\widetilde{\Gamma}_{bc}^{a}(\widetilde{\xi}(t))
\widetilde{\eta}^{b}(t)\widetilde{\eta}^{c}(t).
\]

This shows that, in the induced coordinates associated with $(\widetilde{U},\widetilde{\phi})$, the curve $\zeta$ satisfies the differential system corresponding to the vector field $\mathbf{G}_{\widetilde{U}}$. Therefore, the local vector fields $\mathbf{G}_U$ and $\mathbf{G}_{\widetilde{U}}$ agree on overlaps and define a global smooth vector field $\mathbf{G}\in \mathfrak{X}(TM)$.

Let us now verify that this global vector field has the stated property. Let
$\zeta\colon I\longrightarrow TM$ be a smooth curve and write
$\zeta(t)=\bigl(\gamma(t),\mathbf{V}(t)\bigr)$. First suppose that
$\zeta$ is an integral curve of $\mathbf{G}$. For fixed $t_{0}\in I$, there exist
a chart $(U,\phi)$ of $M$ and an open interval $J\subseteq I$ with
$t_{0}\in J$ such that $\zeta(J)\subseteq \pi^{-1}(U)$. Then
$\zeta|_{J}$ is an integral curve of $\mathbf{G}_U$ and, by what we have proved,
$\mathbf{V}(t)=\dot\gamma(t)$ and $\frac{D\mathbf{V}}{dt}(t)=0$ for every
$t\in J$. Since $t_{0}$ is arbitrary, we conclude that these equalities hold
for every $t\in I$.

Conversely, suppose that $\mathbf{V}(t)=\dot\gamma(t)$ and
$\frac{D\mathbf{V}}{dt}(t)=0$ for every $t\in I$. Given $t_{0}\in I$,
choose again a chart $(U,\phi)$ and an open interval
$J\subseteq I$ with $t_{0}\in J$ such that
$\zeta(J)\subseteq \pi^{-1}(U)$. By the local result already proved, the
restriction $\zeta|_{J}$ is an integral curve of $\mathbf{G}_U$ and
hence of $\mathbf{G}$. Since this holds for every $t_{0}\in I$, it follows
that $\zeta$ is an integral curve of $\mathbf{G}$ on all of $I$.

This proves the characterizing property of $\mathbf{G}$.

For uniqueness, let $\mathbf{H}\in \mathfrak{X}(TM)$ be another smooth vector field
with the same property. Let $(U,\phi)$ be a smooth chart of $M$ and let
$q\in \pi^{-1}(U)$. Consider the integral curve $\zeta$ of $\mathbf{G}$
such that $\zeta(0)=q$. By the characterizing property, $\zeta$ satisfies
$\mathbf{V}(t)=\dot\gamma(t)$ and $\frac{D\mathbf{V}}{dt}(t)=0$ and
is therefore also an integral curve of $\mathbf{H}$. Consequently,
$\mathbf{G}_q=\zeta'(0)=\mathbf{H}_q$. Since $q$ is arbitrary, we conclude
that $\mathbf{G}=\mathbf{H}$.

The last assertion follows immediately: if
$\zeta(t)=\bigl(\gamma(t),\dot\gamma(t)\bigr)$, then the condition
$\mathbf{V}(t)=\dot\gamma(t)$ and $\frac{D\mathbf{V}}{dt}(t)=0$ is equivalent to
$\frac{D\dot\gamma}{dt}(t)=0$, which is precisely the geodesic
equation.
\end{proof}

A geodesic $\gamma\colon I\longrightarrow M$ that cannot be extended to a geodesic on a larger interval is called a \textit{\textbf{maximal geodesic}}. A \textit{\textbf{geodesic segment}} is a geodesic whose domain is a compact interval. The existence of maximal geodesics is established in the following theorem:

\begin{theorem}[Existence, uniqueness, and maximality of geodesics]\label{teo: existencia unicidad geodesicas}\index{existence uniqueness and maximality of geodesics@existence, uniqueness, and maximality of geodesics}
Let $M$ be a smooth manifold without boundary and let $\nabla$ be a connection on $TM$. For each $p\in M$, $v\in T_{p}M$, and $t_{0}\in\mathbb{R}$, there exists a unique maximal geodesic
\[
\gamma_{(p,v,t_{0})}\colon I_{(p,v,t_{0})}\longrightarrow M
\]
such that $t_{0}\in I_{(p,v,t_{0})}$ and $\gamma_{(p,v,t_{0})}(t_{0})=p$, $\dot\gamma_{(p,v,t_{0})}(t_{0})=v$.

Moreover, if $\gamma_{1}\colon I_{1}\longrightarrow M$ and $\gamma_{2}\colon I_{2}\longrightarrow M$ are geodesics such that $\gamma_{1}(t_{0})=\gamma_{2}(t_{0})$ and $\dot\gamma_{1}(t_{0})=\dot\gamma_{2}(t_{0})$, then $\gamma_{1}(t)=\gamma_{2}(t)$ for every $t\in I_{1}\cap I_{2}$.
\end{theorem}

\begin{proof}
Let $\mathbf{G}\in\mathfrak{X}(TM)$ be the geodesic vector field given by Proposition~\ref{prop: campo geodesico}. By Theorem~\ref{teo: fundamental de flujos}, for each $(p,v)\in TM$ there exists a unique maximal integral curve $\zeta_{(p,v,t_{0})}\colon I_{(p,v,t_{0})}\longrightarrow TM$ such that $\zeta_{(p,v,t_{0})}(t_{0})=(p,v)$.

Define $\gamma_{(p,v,t_{0})}:=\pi\circ \zeta_{(p,v,t_{0})}$. Then $\gamma_{(p,v,t_{0})}$ is a smooth curve in $M$ and $\gamma_{(p,v,t_{0})}(t_{0})=\pi(\zeta_{(p,v,t_{0})}(t_{0}))=p$.

Write
$\zeta_{(p,v,t_{0})}(t)
=\bigl(\gamma_{(p,v,t_{0})}(t),\mathbf{V}(t)\bigr)$ for
$t\in I_{(p,v,t_{0})}$, where
$\mathbf{V}(t)\in T_{\gamma_{(p,v,t_{0})}(t)}M$. By
Proposition~\ref{prop: campo geodesico}, we have
$\mathbf{V}(t)=\dot\gamma_{(p,v,t_{0})}(t)$ for every
$t\in I_{(p,v,t_{0})}$, and moreover
\[
\displaystyle\frac{D\dot\gamma_{(p,v,t_{0})}}{dt}(t)=0.
\]
Therefore, $\gamma_{(p,v,t_{0})}$ is a geodesic. Evaluating at $t=t_{0}$
gives $\mathbf{V}(t_{0})=v$, that is,
$\dot\gamma_{(p,v,t_{0})}(t_{0})=v$.

This proves existence.

We now prove uniqueness. Let $\widetilde{\gamma}\colon J\longrightarrow M$ be a geodesic such that $\widetilde{\gamma}(t_{0})=p$ and $\dot{\widetilde{\gamma}}(t_{0})=v$. Define $\widetilde{\zeta}\colon J\longrightarrow TM$ by $\widetilde{\zeta}(t):=\dot{\widetilde{\gamma}}(t)$. Since $\widetilde{\gamma}$ is a geodesic, Proposition~\ref{prop: campo geodesico} implies that $\widetilde{\zeta}$ is an integral curve of $\mathbf{G}$. Moreover, $\widetilde{\zeta}(t_{0})=(p,v)=\zeta_{(p,v,t_{0})}(t_{0})$. Uniqueness of integral curves of the vector field $\mathbf{G}$ gives $\widetilde{\zeta}(t)=\zeta_{(p,v,t_{0})}(t)$ for every $t\in J\cap I_{(p,v,t_{0})}$. Applying the projection $\pi$ gives $\widetilde{\gamma}(t)=\gamma_{(p,v,t_{0})}(t)$ for every $t\in J\cap I_{(p,v,t_{0})}$. This proves uniqueness.

Finally, let us verify maximality. Suppose that $\gamma_{(p,v,t_{0})}$ is not maximal. Then there exists a geodesic $\overline{\gamma}\colon \widetilde{I}\longrightarrow M$ such that $I_{(p,v,t_{0})}\subsetneq \widetilde{I}$ and $\overline{\gamma}(t)=\gamma_{(p,v,t_{0})}(t)$ for every $t\in I_{(p,v,t_{0})}$. Define $\overline{\zeta}\colon \widetilde{I}\longrightarrow TM$ by $\overline{\zeta}(t):=\dot{\overline{\gamma}}(t)$. Since $\overline{\gamma}$ is a geodesic, Proposition~\ref{prop: campo geodesico} implies that $\overline{\zeta}$ is an integral curve of $\mathbf{G}$. Moreover, $\overline{\zeta}(t)=\dot{\overline{\gamma}}(t)=\dot\gamma_{(p,v,t_{0})}(t)=\zeta_{(p,v,t_{0})}(t)$ for every $t\in I_{(p,v,t_{0})}$. Therefore, $\overline{\zeta}$ is a proper extension of the integral curve $\zeta_{(p,v,t_{0})}$, contradicting its maximality. We conclude that $\gamma_{(p,v,t_{0})}$ is maximal.
\end{proof}

From now on, we denote this unique maximal geodesic by $\gamma_{v}$ and assume $t_{0}=0$ to simplify the notation.

\begin{remark}\label{obs: flujo geodesico y accion de G sobre funciones}
By Theorem~\ref{teo: existencia unicidad geodesicas}, there exists a unique maximal flow $\Phi_{t}$ associated with the geodesic vector field $\mathbf{G}$, given by
\[
\Phi_{t}(p,v)=\bigl(\gamma_{v}(t),\dot\gamma_{v}(t)\bigr).
\]

By the definition of the vector field as the derivative of its flow, we have
\[
\mathbf{G}_{(p,v)}
=
\displaystyle\frac{d}{dt}\Phi_{t}(p,v)\Big|_{t=0}.
\]
Therefore, for every $f\in C^{\infty}(TM)$,
\[
\mathbf{G}f(p,v)
=
df_{(p,v)}\bigl(\mathbf{G}_{(p,v)}\bigr)
=
\displaystyle\frac{d}{dt}\Big(f(\Phi_{t}(p,v))\Big)\Big|_{t=0}.
\]

Substituting the expression for the flow gives
\[
\mathbf{G}f(p,v)
=
\displaystyle\frac{d}{dt}\Big(f(\gamma_{v}(t),\dot\gamma_{v}(t))\Big)\Big|_{t=0}.
\]

On the other hand, in induced coordinates $(x^{1},\dots,x^{n},v^{1},\dots,v^{n})$ on $\pi^{-1}(U)$, we have
\[
\mathbf{G}
=
v^{k}\displaystyle\boldsymbol{\partial}_{k}
-
v^{i}v^{j}\Gamma_{ij}^{k}(x)\displaystyle\boldsymbol{\partial}_{v^{k}},
\]
and consequently, for every $f\in C^{\infty}(\pi^{-1}(U))$,
\[
\mathbf{G}f
=
v^{k}\displaystyle\frac{\partial f}{\partial x^{k}}
-
v^{i}v^{j}\Gamma_{ij}^{k}(x)\displaystyle\frac{\partial f}{\partial v^{k}}.
\]

In particular, if $f\in C^{\infty}(M)$, then its lift $\pi^{*}f=f\circ\pi$ satisfies
\[
\mathbf{G}(\pi^{*}f)(p,v)
=
\displaystyle\frac{d}{dt}\Big(f(\gamma_{v}(t))\Big)\Big|_{t=0}
=
df_{p}(v).
\]
\end{remark}

The exponential map compares a neighborhood of the origin in $T_pM$ with a neighborhood of $p$ in the manifold. Before constructing it, we should observe how a geodesic changes when its initial velocity is rescaled.

\begin{lemma}[Rescaling lemma]\label{lema: reescalamiento geodesicas}\index{rescaling lemma!for geodesics}
Let $M$ be a smooth manifold without boundary and let $\nabla$ be a connection on
$TM$. For every $p\in M$, $v\in T_{p}M$, and $c\in\mathbb{R}$, we have
\[
\gamma_{cv}(t)=\gamma_{v}(ct)
\]
for every $t$ for which both sides are defined.
\end{lemma}

\begin{proof}
If $c=0$, both curves are constant and equal to $p$, so the assertion is trivial. Suppose $c\neq 0$.

Let $I\subseteq\mathbb{R}$ be the maximal domain of $\gamma_{v}$ and define
\[
\widetilde{\gamma}_{v}\colon c^{-1}I\longrightarrow M,
\qquad
\widetilde{\gamma}_{v}(t)=\gamma_{v}(ct).
\]
We will show that $\widetilde{\gamma}_{v}$ is a geodesic with initial point $p$ and initial velocity $cv$. By uniqueness and maximality of $\gamma_{cv}$, it will follow that $\widetilde{\gamma}_{v}=\gamma_{cv}$.

Both curves start at $p$, since $\widetilde{\gamma}_{v}(0)=\gamma_{v}(0)=p$. Consider a chart $(U,\phi)$ and suppose that
\[
(\phi\circ \gamma_{v})(t)=(x^{1}(t),\dots,x^{n}(t)).
\]
Then
\[
(\phi\circ \widetilde{\gamma}_{v})(t)=(x^{1}(ct),\dots,x^{n}(ct)).
\]
By the chain rule,
\[
\dot{\widetilde{\gamma}}_{v}^{i}(t)
=
\displaystyle\frac{d}{dt}x^{i}(ct)
=
c\,\dot{x}^{i}(ct).
\]
In particular,
\[
\dot{\widetilde{\gamma}}_{v}(0)=c\,\dot\gamma_{v}(0)=cv.
\]

To see that $\widetilde{\gamma}_{v}$ is a geodesic, we compute its covariant derivative. Using the local expression,
\[
\displaystyle\frac{D\dot{\widetilde{\gamma}}_{v}}{dt}(t)
=
\left(
\displaystyle\frac{d}{dt}\dot{\widetilde{\gamma}}_{v}^{k}(t)
+
\Gamma_{ij}^{k}(\widetilde{\gamma}_{v}(t))
\dot{\widetilde{\gamma}}_{v}^{i}(t)\dot{\widetilde{\gamma}}_{v}^{j}(t)
\right)\boldsymbol{\partial}_{k}.
\]
Substituting the preceding expression for $\dot{\widetilde{\gamma}}_{v}^{i}(t)$ gives
\[
\displaystyle\frac{d}{dt}\dot{\widetilde{\gamma}}_{v}^{k}(t)
=
c^{2}\ddot{x}^{k}(ct),
\]
and
\[
\dot{\widetilde{\gamma}}_{v}^{i}(t)\dot{\widetilde{\gamma}}_{v}^{j}(t)
=
c^{2}\dot{x}^{i}(ct)\dot{x}^{j}(ct).
\]
Therefore,
\[
\displaystyle\frac{D\dot{\widetilde{\gamma}}_{v}}{dt}(t)
=
c^{2}
\left(
\ddot{x}^{k}(ct)
+
\Gamma_{ij}^{k}(\gamma_{v}(ct))
\dot{x}^{i}(ct)\dot{x}^{j}(ct)
\right)\boldsymbol{\partial}_{k}.
\]
Since $\gamma_{v}$ is a geodesic, the term in parentheses is zero, and consequently
\[
\displaystyle\frac{D\dot{\widetilde{\gamma}}_{v}}{dt}(t)=0.
\]

This shows that $\widetilde{\gamma}_{v}$ is a geodesic with initial conditions $(p,cv)$, so $\widetilde{\gamma}_{v}=\gamma_{cv}$.
\end{proof}

\subsection{The exponential map}

The assignment $v\mapsto \gamma_{v}$ defines a map from $TM$ to the set of geodesics of $M$. Moreover, by the rescaling lemma, we may define a function from a subset of $TM$ to $M$ that sends each line through the origin of $T_{p}M$ to a geodesic:

\begin{definition}\label{def:variedades-riemannianas-mapeo-exponencial}\index{exponential map}
Let $M$ be a smooth manifold without boundary and let $\nabla$ be a connection on
$TM$. By Theorem~\ref{teo: existencia unicidad geodesicas}, each
$p\in M$ and $v\in T_{p}M$ determine a unique maximal geodesic
$\gamma_{v}$, defining a map $v\mapsto \gamma_{v}$ from $TM$ to the
set of geodesics in $M$.

Let
\[
\mathcal{E}=\{v\in TM\mid \gamma_{v}\text{ is defined on an interval containing $[0,1]$}\}.
\]
Define the \textbf{exponential map} $\operatorname{exp}\colon \mathcal{E}\longrightarrow M$ by $\operatorname{exp}(v)=\gamma_{v}(1)$. For each $p\in M$, the \textbf{exponential map restricted to $p$}, denoted by $\operatorname{exp}_{p}$, is the restriction of the exponential map to the set $\mathcal{E}_{p}=\mathcal{E}\cap T_{p}M$.
\end{definition}
\begin{figura}[H]
    \includegraphics[width=0.5\linewidth]{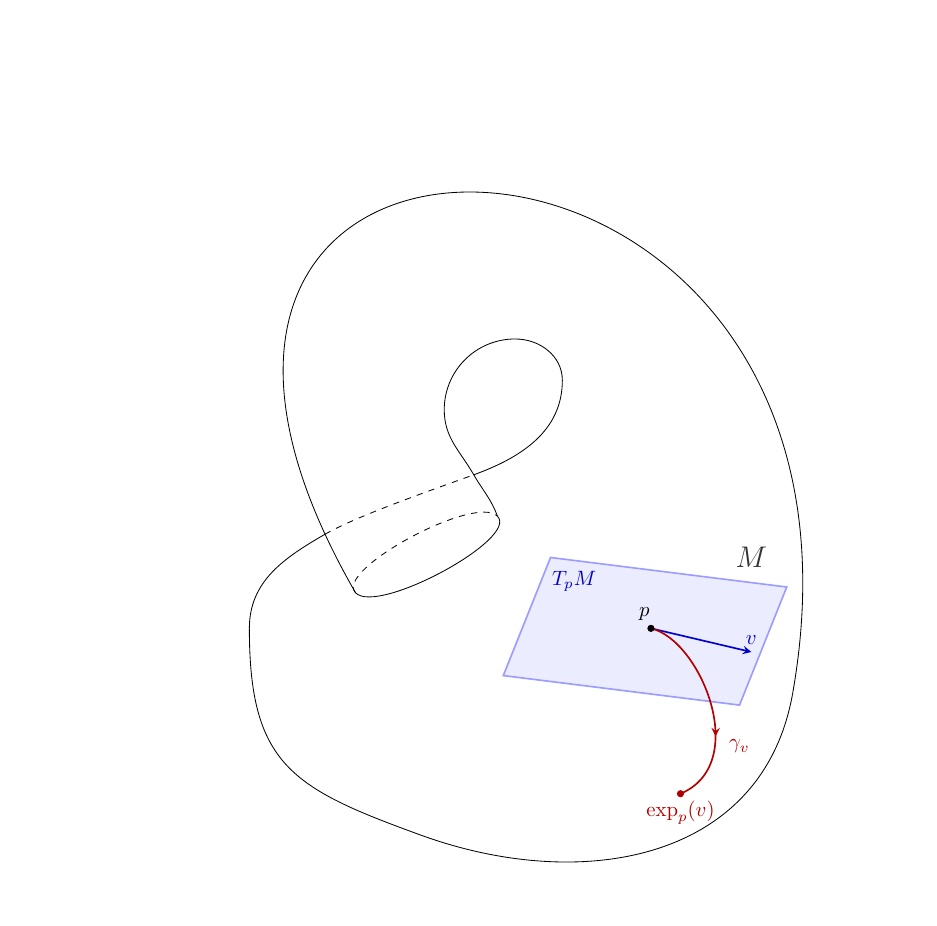}
    \caption{Exponential map on a Riemannian manifold.}
    \label{fig:mapeo-exponencial}
\end{figura}

From this point onward, when $(M,\mathbf g)$ is a Riemannian
manifold, geodesics and the exponential map are understood to be associated
with the Levi--Civita connection of $\mathbf g$, unless otherwise stated.
The preceding results formulated for an arbitrary connection
retain their original meaning.

The following proposition collects the basic properties of the exponential map that will be used to construct normal coordinates and study the local geometry.

\begin{proposition}[Properties of the exponential map]\label{propiedades mapeo exponencial}\index{properties of the exponential map}
Let $(M,\mathbf{g})$ be a Riemannian manifold without boundary and let
$\operatorname{exp}\colon \mathcal{E}\longrightarrow M$ be the exponential
map of its Levi--Civita connection.
\begin{enumerate}
\item[(a)] $\mathcal{E}$ is an open subset of $TM$ containing $0_{p}\in T_{p}M$ for every $p\in M$, and each $\mathcal{E}_{p}\subseteq T_{p}M$ is star-shaped with respect to $0$.
\item[(b)] For each $v\in TM$, the geodesic $\gamma_{v}$ is given by
\[
\gamma_{v}(t)=\operatorname{exp}(tv)
\]
for every $t$ for which both sides are defined.
\item[(c)] $\operatorname{exp}$ is smooth.
\item[(d)] For each $p\in M$, the differential
\[
d(\operatorname{exp}_{p})_{0}\colon T_{0}(T_{p}M)\cong T_{p}M\longrightarrow T_{p}M
\]
is the identity under the usual identification of $T_{0}(T_{p}M)$ with $T_{p}M$.
\end{enumerate}
\end{proposition}

\begin{proof}
By Lemma~\ref{lema: reescalamiento geodesicas}, for every $v\in TM$ and every $t\in\mathbb{R}$ for which both sides are defined, we have
\[
\gamma_{v}(t)=\gamma_{tv}(1)=\operatorname{exp}(tv).
\]
This proves (b).

Moreover, for each $p\in M$, the constant geodesic at $p$ has initial velocity $0_{p}$, so $0_{p}\in\mathcal{E}_{p}$. Now, if $v\in\mathcal{E}_{p}$, then $\gamma_{v}$ is defined on $t=1$. By Lemma~\ref{lema: reescalamiento geodesicas}, for every $t\in[0,1]$ we have
\[
\operatorname{exp}_{p}(tv)=\gamma_{tv}(1)=\gamma_{v}(t),
\]
and in particular $tv\in\mathcal{E}_{p}$. This shows that $\mathcal{E}_{p}$ is star-shaped with respect to $0$.

We now show that $\mathcal{E}$ is open and $\operatorname{exp}$ is smooth. Let $\mathbf{G}$ be the geodesic vector field and let
\[
\Phi\colon \mathcal{D}(\Phi)\longrightarrow TM
\]
be its maximal flow. By Theorem~\ref{teo: fundamental de flujos}, the domain $\mathcal{D}(\Phi)\subseteq \mathbb{R}\times TM$ is open and the map $\Phi$ is smooth. By Remark~\ref{obs: flujo geodesico y accion de G sobre funciones}, we have
\[
\Phi_{t}(v)=\bigl(\gamma_{v}(t),\dot\gamma_{v}(t)\bigr)
\]
whenever both sides are defined.

By the definition of $\mathcal{E}$,
\[
\mathcal{E}
=
\{\,v\in TM \mid \gamma_{v}(1)\text{ is defined}\,\}.
\]
Since the maximal domain of $\gamma_{v}$ contains $0$, it follows that $\gamma_{v}(1)$ is defined if and only if $(1,v)\in\mathcal{D}(\Phi)$. Therefore,
\[
\mathcal{E}
=
\{\,v\in TM \mid (1,v)\in\mathcal{D}(\Phi)\,\}.
\]

If we define $\iota\colon TM\longrightarrow\mathbb{R}\times TM$ by $\iota(v)=(1,v)$, then
\[
\mathcal{E}=\iota^{-1}(\mathcal{D}(\Phi)).
\]
Since $\iota$ is continuous and $\mathcal{D}(\Phi)$ is open, we conclude that $\mathcal{E}$ is open.

Moreover, for every $v\in\mathcal{E}$,
\[
\operatorname{exp}(v)=\gamma_{v}(1)=\pi\bigl(\Phi_{1}(v)\bigr).
\]
Consequently,
\[
\operatorname{exp}=\pi\circ \Phi_{1}\big|_{\mathcal{E}},
\]
and since $\Phi$ and $\pi$ are smooth, it follows that $\operatorname{exp}$ is smooth. This proves (c) and also completes the proof of (a).

Finally, fix $p\in M$ and let $v\in T_{p}M$. Consider the curve $\alpha(t)=tv\in T_{p}M$. By (b),
\[
(\operatorname{exp}_{p}\circ \alpha)(t)=\operatorname{exp}_{p}(tv)=\gamma_{v}(t)
\]
for every $t$ in an open interval around $0$ on which the expression is defined. Differentiating at $t=0$ gives
\[
d(\operatorname{exp}_{p})_{0}(v)
=
\displaystyle\frac{d}{dt}\Big(\operatorname{exp}_{p}(tv)\Big)\Big|_{t=0}
=
\dot\gamma_{v}(0)
=
v.
\]
Under the usual identification of $T_{0}(T_{p}M)$ with $T_{p}M$, this shows that $d(\operatorname{exp}_{p})_{0}$ is the identity. This proves (d).
\end{proof}

\begin{figura}[H]
    \includegraphics[scale=0.65]{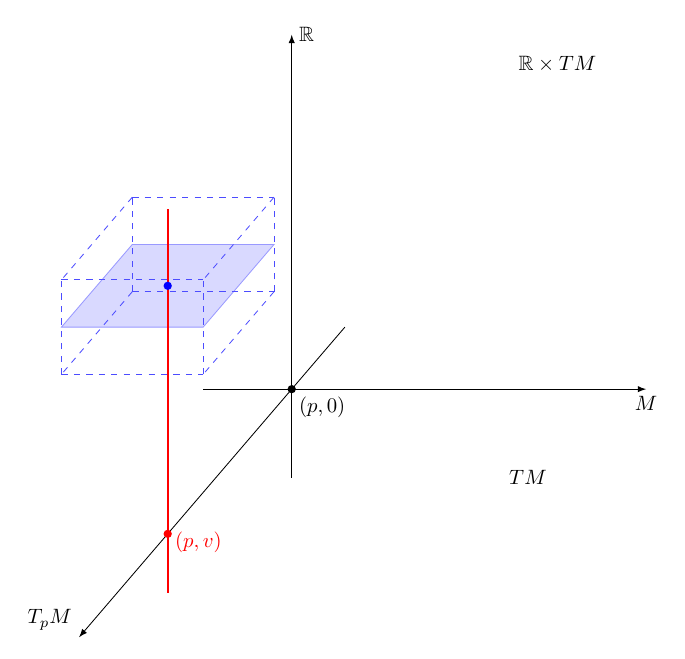}
    \caption{$\mathcal{E}$ is open.}
    \label{fig:E-abierto}
\end{figura}

\subsection{Geodesic normal coordinates}

Among the properties of the exponential map, we proved that $d(\operatorname{exp}_{p})_{0}$ is invertible, so the inverse function theorem for manifolds (Theorem~\ref{teo: funcion inversa variedades}) guarantees the existence of open sets $V\subseteq T_{p}M$ and $U\subseteq M$ such that $0\in V$, $p\in U$, and
\[
\operatorname{exp}_{p}\colon V\longrightarrow U
\]
is a diffeomorphism.

If $V$ is star-shaped with respect to the origin, we say that $U$ is a \textit{\textbf{normal neighborhood of $p$}}. Every orthonormal basis $\{\beta_{1},\dots ,\beta_{n}\}$ of $T_{p}M$ determines a canonical isomorphism $B\colon \mathbb{R}^{n}\longrightarrow T_{p}M$ by $B(x^{1},\dots, x^{n})=\displaystyle\sum_{i=1}^n x^{i}\beta_{i}$. If $U=\operatorname{exp}_{p}(V)$ is a normal neighborhood of $p$, we may combine this isomorphism with the exponential map to obtain a smooth chart $(U,\phi)$ of $M$ containing $p$, where $\phi=B^{-1}\circ (\operatorname{exp}_{p}\restriction_{V})^{-1}\colon U\longrightarrow \mathbb{R}^{n}$. These coordinates are called \textit{\textbf{normal coordinates centered at $p$}}. They are unique in the following sense:

\begin{proposition}\label{unicidad coordenadas normales}
Let $(M,\mathbf{g})$ be a Riemannian manifold of dimension $n$ without boundary, let $p\in M$, and let $U$ be a normal neighborhood of $p$. For each normal coordinate chart $(U,\phi)$ centered at $p$, the coordinate basis $\{\boldsymbol{\partial}_{1}|_{p},\dots,\boldsymbol{\partial}_{n}|_{p}\}\subseteq T_{p}M$ is orthonormal. Conversely, for each orthonormal basis $\{\beta_{1},\dots,\beta_{n}\}$ of $T_{p}M$, there exists a unique normal coordinate chart $(x^{i})$ on $U$ such that $\boldsymbol{\partial}_{i}|_{p}=\beta_{i}$ for every $i\in\{1,\dots,n\}$. Moreover, any two normal coordinate charts $(x^{i})$ and $(\widetilde{x}^{i})$ are related by
\[
\widetilde{x}^{j}=A_{i}^{j}x^{i}
\]
for some matrix $(A_{i}^{j})\in O(n)$.
\end{proposition}

\begin{proof}
Let $(U,\phi)$ be a normal coordinate chart centered at $p$, with coordinate functions $\phi=(x^{1},\dots,x^{n})$. By definition,
\[
\phi=B^{-1}\circ \operatorname{exp}_{p}^{-1},
\]
where $B\colon \mathbb{R}^{n}\longrightarrow T_{p}M$ is an isomorphism determined by an orthonormal basis $\{b_{1},\dots,b_{n}\}$ of $T_{p}M$. Then
\[
\phi^{-1}=\operatorname{exp}_{p}\circ B,
\]
and therefore
\[
d(\phi^{-1})_{0}=d(\operatorname{exp}_{p})_{0}\circ dB_{0}=B,
\]
since $d(\operatorname{exp}_{p})_{0}=\operatorname{Id}_{T_{p}M}$ and $B$ is linear. Consequently,
\[
\boldsymbol{\partial}_{i}|_{p}=d(\phi^{-1})_{0}(\boldsymbol{\partial}_{i}|_{0})=B(\boldsymbol{\partial}_{i}|_{0})=b_{i},
\]
which shows that the basis $\{\boldsymbol{\partial}_{i}|_{p}\}$ is orthonormal.

Conversely, an orthonormal basis $\{\beta_{1},\dots,\beta_{n}\}$ of $T_{p}M$ determines an isomorphism $B\colon \mathbb{R}^{n}\longrightarrow T_{p}M$ and hence a normal coordinate chart $\phi=B^{-1}\circ \operatorname{exp}_{p}^{-1}$ satisfying $\boldsymbol{\partial}_{i}|_{p}=\beta_{i}$.

If $\widetilde{\phi}=\widetilde{B}^{-1}\circ\operatorname{exp}_{p}^{-1}$ is another normal coordinate chart centered at $p$, then the change of coordinates is given by
\[
\widetilde{\phi}\circ \phi^{-1}
=
\widetilde{B}^{-1}\circ B,
\]
which is a linear isometry of $\mathbb{R}^{n}$. Therefore, there exists a constant matrix $(A_{i}^{j})\in O(n)$ such that
\[
\widetilde{x}^{j}=A_{i}^{j}x^{i}.
\]
This proves uniqueness.
\end{proof}

Other properties of normal coordinates include the following:

\begin{proposition}[Properties of normal coordinates]\label{propiedades de las coordenadas normales}\index{properties of normal coordinates}
Let $(M,\mathbf{g})$ be a Riemannian manifold without boundary and let $(U,\phi)$ be a normal coordinate chart centered at $p\in M$.
\begin{enumerate}
 \item The coordinates of $p$ are $(0,\dots,0)$.
 \item The metric components at $p$ are $g_{ij}(p)=\delta_{ij}$.
 \item For $v=v^{i}\boldsymbol{\partial}_{i}|_{p}\in T_{p}M$, the geodesic $\gamma_{v}$ starting at $p$ with $\dot\gamma(0)=v$ is represented in normal coordinates by
 \[
 (\phi \circ \gamma_{v})(t)=(t v^{1},\dots,t v^{n})
 \]
 if $t\in I$, where $I$ is an interval such that $0\in I$ and $\gamma_{v}(I)\subseteq U$.
 \item The Christoffel symbols in these coordinates vanish at $p$.
 \item The first derivatives of the metric components vanish at $p$:
 \[
 \partial_{k}g_{ij}(p)=0
 \quad\text{for all } i,j,k.
 \]
\end{enumerate}
\end{proposition}

\begin{proof}
Part (a) follows directly from the definition of normal coordinates, since $\phi(p)=0$.

For (b), Proposition~\ref{unicidad coordenadas normales} shows that the coordinate basis $\{\boldsymbol{\partial}_{i}|_{p}\}$ is orthonormal, giving $g_{ij}(p)=\delta_{ij}$.

Part (c) follows from Proposition~\ref{propiedades mapeo exponencial}(b): since $\gamma_{v}(t)=\operatorname{exp}_{p}(tv)$, composing with $\phi=B^{-1}\circ \operatorname{exp}_{p}^{-1}$ gives
\[
(\phi\circ\gamma_{v})(t)=B^{-1}(tv)=(t v^{1},\dots,t v^{n}).
\]

To prove (d), let $v=v^{i}\boldsymbol{\partial}_{i}|_{p}\in T_{p}M$. By (c), writing $x^{i}(t)=(\phi\circ \gamma_{v})^{i}(t)$ gives $x^{i}(t)=t v^{i}$, so $\dot{x}^{i}(t)=v^{i}$ and $\ddot{x}^{k}(t)=0$. Substituting into the equation for the geodesic \eqref{eq: ecuacion de la geodesica}, we obtain
\[
\Gamma_{ij}^{k}(x(t))\,v^{i}v^{j}=0.
\]
Evaluating at $t=0$, that is, at $p$, we conclude that
\[
\Gamma_{ij}^{k}(p)\,v^{i}v^{j}=0
\]
for every $k$ and every $v\in T_{p}M$.

The Levi--Civita connection is torsion-free; therefore, in a
coordinate frame,
\[
 \Gamma_{ab}^{k}=\Gamma_{ba}^{k}.
\]
Taking $v=\boldsymbol{\partial}_{a}|_p$ gives
$\Gamma_{aa}^{k}(p)=0$. If $a\neq b$, taking
$v=\boldsymbol{\partial}_{a}|_p+\boldsymbol{\partial}_{b}|_p$ gives
\[
 0=\Gamma_{aa}^{k}(p)+\Gamma_{ab}^{k}(p)
   +\Gamma_{ba}^{k}(p)+\Gamma_{bb}^{k}(p)
  =2\Gamma_{ab}^{k}(p).
\]
Thus, $\Gamma_{ab}^{k}(p)=0$ for all $a,b,k$, proving (d).

Finally, the Levi--Civita connection is compatible with $\mathbf g$.
Corollary~\ref{compatibilidad levi civita metrica en coordenadas},
applied to $\mathbf E_k=\boldsymbol{\partial}_k$ and evaluated at $p$, gives
\[
 \partial_k g_{ij}(p)
 =\Gamma_{ki}^{a}(p)g_{aj}(p)
  +\Gamma_{kj}^{a}(p)g_{ia}(p)=0,
\]
by (d). This proves (e).
\end{proof}

Continuity of the metric allows us to make the identity $g_{ij}(p)=\delta_{ij}$ quantitative on a sufficiently small ball. The following result specifies the type of normal coordinates we will use in approximation arguments.

\begin{proposition}[Normal coordinate balls with Euclidean control]
\label{prop:bolas-coordenadas-normales-control-euclidiano}
\index{normal coordinate ball!Euclidean control}
Let $(M,\mathbf{g})$ be a Riemannian manifold with or without boundary, let $p\in\operatorname{Int}(M)$, and let $0<\theta<1$. There exist numbers $0<r<R$ and a linear isometry $B_p\colon \mathbb R^n\longrightarrow(T_pM,\mathbf{g}_p)$ such that $\exp_p$ is a diffeomorphism from $B_{\mathbf{g}_p}(0,R)$ onto an open set $U'_p\subseteq\operatorname{Int}(M)$ and, setting
\[
U_p:=\exp_p\bigl(B_{\mathbf{g}_p}(0,r)\bigr),
\qquad
\kappa_p:=B_p^{-1}\circ\exp_p^{-1}\colon U_p\longrightarrow B_{\mathrm{euc}}(0,r),
\]
the pair $(U_p,\kappa_p)$ is a regular normal coordinate ball centered at $p$, with $\overline U_p\subseteq U'_p$, and satisfies
\[
(1-\theta)\|w\|
\leq
\left|d(\kappa_p^{-1})_z w\right|_{\mathbf{g}}
\leq
(1+\theta)\|w\|
\]
for all $z\in B_{\mathrm{euc}}(0,r)$ and $w\in\mathbb R^n$.
Here $\exp_p$ is the exponential map of the Riemannian manifold
$(\operatorname{Int}M,\mathbf g)$.
\end{proposition}

\begin{proof}
Since $d(\exp_p)_0=\operatorname{Id}_{T_pM}$, the inverse function theorem gives an open neighborhood $W\subseteq T_pM$ of $0$ on which $\exp_p$ is a diffeomorphism. Since $p$ is an interior point, we may shrink $W$ so that its image is contained in $\operatorname{Int}(M)$. Choose $R>0$ so that $\overline B_{\mathbf{g}_p}(0,R)\subseteq W$ and fix a linear isometry $B_p\colon \mathbb R^n\longrightarrow(T_pM,\mathbf{g}_p)$. The normal chart
\[
\kappa'_p:=B_p^{-1}\circ\exp_p^{-1}:
U'_p:=\exp_p\bigl(B_{\mathbf{g}_p}(0,R)\bigr)
\longrightarrow B_{\mathrm{euc}}(0,R)
\]
is centered at $p$. If $G(z)$ denotes the matrix of the metric in these coordinates, Proposition~\ref{propiedades de las coordenadas normales} implies $G(0)=I_n$. Since $z\mapsto G(z)$ is also continuous in the operator norm, there exists $0<r<R$ such that
\[
(1-\theta)^2\|w\|^2
\leq
\langle G(z)w,w\rangle_{\mathbb R^n}
\leq
(1+\theta)^2\|w\|^2
\]
for every $z\in\overline B_{\mathrm{euc}}(0,r)$ and every $w\in\mathbb R^n$. Indeed, it suffices to make $\|G(z)-I_n\|_{\operatorname{op}}$ smaller than the minimum of $1-(1-\theta)^2$ and $(1+\theta)^2-1$.

The restriction of $\kappa'_p$ to $U_p:=\exp_p(B_{\mathbf{g}_p}(0,r))$ is the chart $\kappa_p$ in the statement. Moreover, $\overline U_p=\exp_p(\overline B_{\mathbf{g}_p}(0,r))$ is compact and contained in $U'_p$, so $U_p$ is a regular normal coordinate ball in the sense of Definition~\ref{def:nociones-fundamentales-norma-variedad-suave-label-bola}. Finally,
\[
\left|d(\kappa_p^{-1})_z w\right|_{\mathbf{g}}^2
=
\langle G(z)w,w\rangle_{\mathbb R^n},
\]
and taking square roots gives the required inequalities.
\end{proof}

\section{The Riemannian distance function}

The metric first acts infinitesimally, assigning a length to each tangent vector. Integrating these lengths along curves gives a global notion of distance between points. This step is more than a formality: the resulting distance reflects both the local geometry and the ways in which curves can travel through the manifold.
\begin{definition}\label{def:variedades-riemannianas-segmento-de-curva}\index{curve segment}
 Let $M$ be a Riemannian manifold with or without boundary.
 \begin{enumerate}[label=(\alph*)]
 \item A \textbf{curve segment} is a smooth curve whose domain is a compact interval.
 \item A smooth curve $\gamma\colon I\longrightarrow M$ is called a \textbf{regular curve} if $\gamma'(t)\neq 0$ for every $t\in I$ (which implies that $\gamma$ is an immersion).
 \item A continuous curve $\gamma\colon [a,b]\longrightarrow M$ is said to be \textbf{piecewise regular} if there exists a partition $\{t_{0},\dots,t_{k}\}$ of $[a,b]$ such that $\gamma|_{[t_{i-1},t_{i}]}$ is a regular (smooth) curve segment for each $i\in\{1,\dots,k\}$. We also call these \textbf{admissible curves}, and partitions satisfying this condition are called \textbf{admissible partitions for $\gamma$}.
 We also regard constant curves on any compact interval, including the degenerate interval $\{a\}$, as admissible. Apart from this constant case, an admissible curve is always piecewise regular.
 \item If $\gamma\colon I\longrightarrow M$ is a smooth curve, a \textbf{reparametrization of $\gamma$} is a curve of the form $\widetilde{\gamma}=\gamma\circ\phi\colon I'\longrightarrow M$, where $I'$ is an interval and $\phi\colon I'\longrightarrow I$ is a diffeomorphism. If $\phi$ is strictly increasing, we say that $\widetilde{\gamma}$ is a \textbf{forward parametrization} or an \textbf{orientation-preserving parametrization}. If $\phi$ is strictly decreasing, we say that $\widetilde{\gamma}$ is a \textbf{backward parametrization} or an \textbf{orientation-reversing parametrization}.
 \item For an admissible curve $\gamma\colon [a,b]\longrightarrow M$, a \textbf{reparametrization of $\gamma$} is a curve of the form $\widetilde{\gamma}=\gamma\circ \phi$, where $\phi\colon [c,d]\longrightarrow [a,b]$ is a homeomorphism for which there exists a partition $\{c_{0},\dots,c_{k}\}$ of $[c,d]$ such that the restriction $\phi|_{[c_{i-1},c_{i}]}$ is a diffeomorphism onto its image for each $i\in \{1,\dots,k\}$.
 \item If $\gamma\colon [a,b]\longrightarrow M$ is an admissible curve, we define the \textbf{length of $\gamma$} to be \[L_{\mathbf{g}}(\gamma):=\int_{a}^{b}|\dot\gamma(t)|_{\mathbf{g}_{\gamma(t)}}\,dt.\]

 \end{enumerate}
\end{definition}

When concatenating admissible curves, we omit constant factors, which have zero length. The resulting curve is again admissible.

The same integral defines length for any continuous piecewise smooth curve, even if its velocity vanishes. In the definition of distance we continue to take the infimum over the admissible curves specified above.

\begin{proposition}[Arc length and reparametrization]
\label{prop:longitud-arco-reparametrizacion}
Let \((M,\mathbf g)\) be a Riemannian manifold with or without boundary,
and let \(c\colon[a,b]\longrightarrow M\) be a nonconstant admissible
curve. If \(a=t_0<\cdots<t_N=b\) is an admissible partition, then
\[
 L_{\mathbf g}(c)
 =\displaystyle\sum_{j=1}^{N}
   \int_{t_{j-1}}^{t_j}|\dot c(t)|_{\mathbf g_{c(t)}}\,dt.
\]
This value is independent of the partition and is preserved by
admissible reparametrizations. Moreover, if
\[
 s(t):=\int_a^t|\dot c(u)|_{\mathbf g_{c(u)}}\,du,
 \qquad \ell:=L_{\mathbf g}(c),
\]
then \(s\colon[a,b]\longrightarrow[0,\ell]\) is an increasing
homeomorphism, smooth with smooth inverse on each piece of the
partition, and \(\widetilde c:=c\circ s^{-1}\) has unit speed on every
smooth piece.
\end{proposition}
\begin{proof}
Refining a partition merely splits an integral into integrals over
subintervals; any two partitions have a common refinement. This proves
independence. If \(\phi\colon[c_0,d_0]\longrightarrow[a,b]\) is an
admissible reparametrization, refine its partition to include the
preimages of the \(t_j\). On each piece, the chain rule gives
\[
 |(c\circ\phi)'(u)|_{\mathbf g_{c(\phi(u))}}
 =|\dot c(\phi(u))|_{\mathbf g_{c(\phi(u))}}\,|\phi'(u)|.
\]
The substitution \(t=\phi(u)\), using the absolute value of the
Jacobian also when \(\phi\) is decreasing, shows that corresponding
pieces have the same length. Summing proves invariance.

On each piece, \(s'(t)=|\dot c(t)|_{\mathbf g_{c(t)}}>0\).
Thus \(s\) is continuous and strictly increasing on \([a,b]\), with
endpoint values \(0\) and \(\ell\). It is a homeomorphism onto
\([0,\ell]\), and the one-dimensional inverse function theorem gives
a smooth inverse on each piece, up to its endpoints by local
extensions. Finally,
\[
 |\widetilde c'(\sigma)|_{\mathbf g_{\widetilde c(\sigma)}}
 =\frac{|\dot c(t)|_{\mathbf g_{c(t)}}}{s'(t)}=1,
 \qquad t=s^{-1}(\sigma),
\]
in the interiors of these pieces.
\end{proof}
We may define the distance between two points on a connected Riemannian manifold as follows:
\begin{definition}\label{def:variedades-riemannianas-distancia-riemanniana-de-p-a-q}\index{Riemannian distance from p to q}
 Let $(M,\mathbf{g})$ be a connected Riemannian manifold with or without boundary. For each pair of points $p,q\in M$, define the \textbf{Riemannian distance from $p$ to $q$}, denoted by $d_{\mathbf{g}}(p,q)$, to be \[d_{\mathbf{g}}(p,q):=\displaystyle\inf\{L_{\mathbf{g}}(\gamma)\mid \gamma \text{ is an admissible curve from $p$ to $q$}\}.\]
\end{definition}
The infimum is finite on each connected component by the following result.
\begin{proposition}\label{prop:variedades-riemannianas-variedad-suave-conexa-frontera-dos-puntos}
 If $M$ is a connected smooth manifold (with or without boundary), then any two points of $M$ can be joined by an admissible curve.
\end{proposition}
\begin{proof}
Fix $p\in M$ and let $A$ be the set of points that can be joined to
$p$ by an admissible curve. The set $A$ is nonempty. If
$q\in A$, choose a chart around $q$ whose image contains a
Euclidean ball, or a relative ball in $\mathbb H^n$, centered at the
image of $q$. After shrinking it, this set is convex. Straight
segments in the chart, transported to $M$, join $q$ to all the
points in a neighborhood of it; concatenating them with a curve from $p$ to $q$
gives an admissible curve. Therefore, $A$ is open. The same
argument shows that $M\setminus A$ is open: if one of its points
had a point of $A$ in this neighborhood, the local segment would
also join it to $p$. Since $M$ is connected, it follows that $A=M$.
\end{proof}
Every Riemannian manifold is metrizable by its Riemannian distance.
\begin{theorem}\label{variedad riemanniana espacio métrico}\index{Riemannian distance!metric structure of the manifold}
 Let $(M,\mathbf{g})$ be a connected Riemannian manifold with or without boundary. With the distance function $d_{\mathbf{g}}$, $M$ is a metric space whose metric topology agrees with its original topology as a smooth manifold.
\end{theorem}
\begin{proof}
The preceding proposition shows that $d_{\mathbf{g}}(p,q)$ is finite. Nonnegativity
is immediate. Reversing the parameter of a curve proves
that $d_{\mathbf{g}}(p,q)=d_{\mathbf{g}}(q,p)$, and concatenating admissible curves, then
taking infima, gives
\[
 d_{\mathbf{g}}(p,r)\leq d_{\mathbf{g}}(p,q)+d_{\mathbf{g}}(q,r).
\]

It remains to prove that distinct points have positive distance and to compare
the topologies. Fix $p\in M$. In a chart $(U,x)$ around $p$,
choose $r>0$ so that the relative ball
\[
 Q:=B_{\mathrm{euc}}(x(p),2r)\cap\mathbb H^n
\]
is contained in $x(U)$ and has closure contained in $x(U)$; in an
interior chart, omit the intersection with $\mathbb H^n$. Continuity
and positive definiteness of the metric give constants
$0<c\leq C$ such that
\begin{equation}
\label{eq:comparacion-local-distancia-coordenadas}
 c|v|\leq
 \left|d(x^{-1})_zv\right|_{\mathbf{g}}
 \leq C|v|
\end{equation}
for $z\in\overline Q$ and $v\in\mathbb R^n$ tangent to the half-space.
Set
\[
 V:=x^{-1}\bigl(B_{\mathrm{euc}}(x(p),r)\cap\mathbb H^n\bigr).
\]

If $q\in V$, the straight segment from $x(p)$ to $x(q)$ stays in the relative
ball, and the second inequality in
\eqref{eq:comparacion-local-distancia-coordenadas} shows that
\[
 d_{\mathbf{g}}(p,q)\leq C|x(q)-x(p)|.
\]
On the other hand, every admissible curve from $p$ to a point outside $V$
must first reach the coordinate sphere of radius $r$. On the portion
before this first exit point, the first inequality in
\eqref{eq:comparacion-local-distancia-coordenadas} gives a length of at least
$cr$. If $q\in V$ and a curve from $p$ to $q$ does not leave $V$, the same
inequality and the Euclidean inequality between length and displacement
give a length of at least $c|x(q)-x(p)|$; if it leaves $V$, its length is
at least $cr$ and the same bound remains valid. Consequently,
\begin{equation}
\label{eq:bilipschitz-local-distancia-riemanniana}
 c|x(q)-x(p)|\leq d_{\mathbf{g}}(p,q)
 \leq C|x(q)-x(p)|,
 \qquad q\in V.
\end{equation}

These bounds first prove that $d_{\mathbf{g}}(p,q)>0$ when $p\ne q$: if
$q\notin V$, use the bound $cr$, and if $q\in V$, use the first inequality
in \eqref{eq:bilipschitz-local-distancia-riemanniana}. Therefore, $d_{\mathbf{g}}$ is
a metric. The second inequality shows that convergence in the
manifold topology implies convergence with respect to $d_{\mathbf{g}}$ near each
point; the first, together with the barrier $cr$ for leaving $V$, gives the
converse. The two topologies agree.
\end{proof}

\begin{example}[A function defined by intrinsic data]
\label{ej:distancia-funcion-intrinseca}
Let $(M,\mathbf g)$ be a connected Riemannian manifold with or without
boundary, and fix $q\in M$. The assignment
\[
 r_q\colon M\longrightarrow[0,\infty),\qquad r_q(p):=d_{\mathbf g}(q,p)
\]
defines a function without choosing a Euclidean embedding or coordinates.
Its data are the metric and the chosen point. By the triangle inequality,
\[
 d_{\mathbf g}(q,p)\leq d_{\mathbf g}(q,p')+d_{\mathbf g}(p',p),
 \qquad
 d_{\mathbf g}(q,p')\leq d_{\mathbf g}(q,p)+d_{\mathbf g}(p,p').
\]
Subtracting in both directions gives
\[
 |r_q(p)-r_q(p')|\leq d_{\mathbf g}(p,p').
\]
Thus $r_q$ is $1$-Lipschitz and, by the preceding theorem, is continuous
in the manifold topology. Proposition~\ref{prop:regularidad-escalar-en-cartas}
ensures that its representations in every chart are continuous.

A geometric definition does not imply smoothness. On
$(\mathbb R^n,\mathbf g_{\mathrm{euc}})$, with $n\geq1$ and $q=0$,
$r_0(x)=|x|$: every path has length at least $|x|$, and the straight
segment attains that length. Along the line $x=te_1$,
$r_0(te_1)=|t|$, which is not differentiable at $t=0$.
By contrast, $r_0^2(x)=\displaystyle\sum_{i=1}^n(x^i)^2$ is smooth.
The choice of intrinsic data and the regularity of the function are
distinct matters; both must be specified when formulating a problem.
\end{example}

\begin{definition}\label{def:variedades-riemannianas-curva-minimizante}\index{minimizing curve}
Let \((M,\mathbf g)\) be a Riemannian manifold without boundary. An
admissible curve \(c\colon[a,b]\longrightarrow M\), with endpoints in
the same connected component, is \textbf{minimizing} if
\[
 L_{\mathbf g}(c)=d_{\mathbf g}(c(a),c(b)).
\]
Equivalently, its length does not exceed that of any admissible curve
with the same endpoints; a competing curve may have a different
parameter interval.
\end{definition}

\begin{definition}[Variation and variational field]
\label{def:variacion-campo-variacional-jacobi}
\index{variation of curves}\index{variational field}
Let $(M,\mathbf{g})$ be a Riemannian manifold without boundary. A
\textbf{smooth variation} of a curve
$\gamma\colon I\to M$ is a smooth map
\[
 \alpha\colon(-\varepsilon,\varepsilon)\times I\longrightarrow M,
 \qquad \alpha(0,t)=\gamma(t).
\]
We write
\[
 \mathbf{S}:=d\alpha(\boldsymbol{\partial}_s),\qquad \mathbf{T}:=d\alpha(\boldsymbol{\partial}_t),
\]
and call $\mathbf{J}(t):=\mathbf{S}(0,t)$ the variational field along $\gamma$.
The variation is \textbf{through geodesics} if
$t\mapsto\alpha(s,t)$ is a geodesic for each $s$.
\end{definition}

Covariant differentiation along curves can also be applied to each
variable of a variation. For a smooth section \(\mathbf Z\) of
\(\alpha^*TM\), write \(D_s\mathbf Z\) and \(D_t\mathbf Z\) for these
derivatives, keeping the other variable fixed. In a chart on \(M\),
with \(n=\dim M\), their expressions are
\[
\begin{aligned}
 D_s\mathbf Z
 &=\displaystyle\sum_{k=1}^{n}
 \left(\partial_sZ^k+
 \displaystyle\sum_{i,j=1}^{n}
 \Gamma^k_{ij}(\alpha)S^iZ^j\right)
 \boldsymbol{\partial}_k\big|_\alpha,\\
 D_t\mathbf Z
 &=\displaystyle\sum_{k=1}^{n}
 \left(\partial_tZ^k+
 \displaystyle\sum_{i,j=1}^{n}
 \Gamma^k_{ij}(\alpha)T^iZ^j\right)
 \boldsymbol{\partial}_k\big|_\alpha.
\end{aligned}
\]
These are the formulas of
Proposition~\ref{prop: expresion local derivada covariante curvas}
applied to the curves in each family; in particular, they define smooth
fields along \(\alpha\).

\begin{lemma}[Commutation of the velocities of a variation]
\label{lem:conmutacion-velocidades-variacion}
Let \(\alpha\) be a smooth variation in a Riemannian manifold without
boundary. If
\(\mathbf S=d\alpha(\boldsymbol{\partial}_s)\) and
\(\mathbf T=d\alpha(\boldsymbol{\partial}_t)\), then
\[
 D_s\mathbf T=D_t\mathbf S.
\]
\end{lemma}
\begin{proof}
In a chart, write
\(S^i=\partial_s\alpha^i\) and \(T^i=\partial_t\alpha^i\),
for \(1\leq i\leq n\). For every \(1\leq k\leq n\),
\[
\begin{aligned}
 (D_s\mathbf T)^k-(D_t\mathbf S)^k
 &=
 \partial_s\partial_t\alpha^k-\partial_t\partial_s\alpha^k\\
 &\quad+
 \displaystyle\sum_{i,j=1}^{n}
 \Gamma^k_{ij}(\alpha)(S^iT^j-T^iS^j)=0.
\end{aligned}
\]
Mixed derivatives commute and
\(\Gamma^k_{ij}=\Gamma^k_{ji}\), since the Levi--Civita connection has
zero torsion. The calculation applies in every chart and proves the
identity throughout the parameter domain.
\end{proof}

\begin{lemma}[Gauss lemma]\label{lem:gauss-exponencial}\index{Gauss lemma}
Let \((M,\mathbf g)\) be a Riemannian manifold without boundary, let
\(p\in M\), and let \(v,w\in T_pM\), with \(v\in\mathcal E_p\).
Identifying \(T_v(T_pM)\) with \(T_pM\), we have
\[
 \mathbf g_{\exp_p(v)}
 \bigl(d(\exp_p)_v(v),d(\exp_p)_v(w)\bigr)
 =\mathbf g_p(v,w).
\]
\end{lemma}
\begin{proof}
Since \(\mathcal E_p\) is open and star-shaped about zero, there is
\(\delta>0\) such that
\[
 \alpha\colon(-\delta,\delta)\times[0,1]\longrightarrow M,
 \qquad \alpha(s,t):=\exp_p\bigl(t(v+sw)\bigr)
\]
is well defined and smooth. For each \(s\), the curve
\(\alpha_s(t)=\alpha(s,t)\) is the geodesic with initial velocity
\(v+sw\). Set
\(\mathbf S=d\alpha(\boldsymbol{\partial}_s)\) and
\(\mathbf T=d\alpha(\boldsymbol{\partial}_t)\).
The geodesic equation and metric compatibility give
\[
 \partial_t|\mathbf T(s,t)|_{\mathbf g_{\alpha(s,t)}}^2
 =2\mathbf g_{\alpha(s,t)}(D_t\mathbf T,\mathbf T)=0.
\]
Consequently,
\(|\mathbf T(s,t)|_{\mathbf g_{\alpha(s,t)}}^2=|v+sw|_{\mathbf g_p}^2\).
Using Lemma~\ref{lem:conmutacion-velocidades-variacion}, we obtain
\[
\begin{aligned}
 \partial_t\bigl(\mathbf g_\alpha(\mathbf S,\mathbf T)\bigr)
 &=\mathbf g_\alpha(D_t\mathbf S,\mathbf T)
   +\mathbf g_\alpha(\mathbf S,D_t\mathbf T)\\
 &=\mathbf g_\alpha(D_s\mathbf T,\mathbf T)
 =\frac12\partial_s|v+sw|_{\mathbf g_p}^2.
\end{aligned}
\]
Restricting to \(s=0\), the last member is \(\mathbf g_p(v,w)\).
Moreover, \(\alpha(s,0)=p\), so \(\mathbf S(0,0)=0\).
Integrating with respect to \(t\) over \([0,1]\), we find
\[
 \mathbf g_{\alpha(0,1)}(\mathbf S(0,1),\mathbf T(0,1))
 =\int_0^1\mathbf g_p(v,w)\,dt
 =\mathbf g_p(v,w).
\]
The chain rule identifies
\(\mathbf S(0,1)=d(\exp_p)_v(w)\) and
\(\mathbf T(0,1)=d(\exp_p)_v(v)\), as required.
\end{proof}

\begin{proposition}[Radial minimization in a normal neighborhood]
\label{prop:minimizacion-radial-vecindad-normal}
Let \((M,\mathbf g)\) be a Riemannian manifold without boundary and let
\(p\in M\). There is \(\varepsilon>0\) such that \(\exp_p\) is a
diffeomorphism from \(B_{\mathbf g_p}(0,\varepsilon)\) onto an open
set \(U\). For every \(\varepsilon\) with this property and every
\(v\in B_{\mathbf g_p}(0,\varepsilon)\), the geodesic
\[
 \gamma_v\colon[0,1]\longrightarrow M,
 \qquad \gamma_v(t)=\exp_p(tv),
\]
satisfies
\[
 L_{\mathbf g}(\gamma_v)
 =d_{\mathbf g}(p,\exp_p(v))=|v|_{\mathbf g_p}.
\]
If \(v\neq0\), it is the unique admissible minimizing curve between
these points, up to orientation-preserving reparametrizations.
If \(v=0\), every curve of minimum length is constant.
Moreover, \(U=B_{\mathbf g}(p,\varepsilon)\).
\end{proposition}
\begin{proof}
The existence of \(\varepsilon\) follows from
\(d(\exp_p)_0=\operatorname{Id}_{T_pM}\) and the inverse function
theorem. Fix any radius satisfying the statement. Define
\[
 \rho\colon U\longrightarrow[0,\varepsilon),
 \qquad \rho(q):=|\exp_p^{-1}(q)|_{\mathbf g_p}.
\]
This function is continuous on \(U\) and smooth on
\(U\setminus\{p\}\). For \(q=\exp_p(z)\neq p\), define the radial
unit vector field by
\[
 \mathbf R(q):=\frac{d(\exp_p)_z(z)}{|z|_{\mathbf g_p}}.
\]
The chain rule and Lemma~\ref{lem:gauss-exponencial} give, for every
\(w\in T_pM\),
\[
\begin{aligned}
 d\rho_q\bigl(d(\exp_p)_z(w)\bigr)
 &=\frac{\mathbf g_p(z,w)}{|z|_{\mathbf g_p}}\\
 &=\mathbf g_q\bigl(\mathbf R(q),d(\exp_p)_z(w)\bigr).
\end{aligned}
\]
The same lemma, with \(w=z\), gives
\(|\mathbf R(q)|_{\mathbf g_q}=1\).
Since \(d(\exp_p)_z\) is surjective, we conclude that
\[
 d\rho_q(\xi)=\mathbf g_q(\mathbf R(q),\xi),
 \qquad |d\rho_q(\xi)|\leq|\xi|_{\mathbf g_q},
 \qquad \xi\in T_qM.
\]

Let \(c\colon[a,b]\longrightarrow U\) be an admissible curve from
\(p\) to \(q=\exp_p(v)\), and set \(r=\rho\circ c\).
We first justify the regularity needed for integration, since \(\rho\)
is not smooth at \(p\). On each regular piece of \(c\), the curve
\(z=\exp_p^{-1}\circ c\) is regular. Its zeros are isolated and form
a closed subset of a compact interval, so there are only finitely many.
Near a zero \(t_0\), the identity
\[
 z(t)=(t-t_0)h(t),\qquad
 h(t):=\int_0^1 z'\bigl(t_0+u(t-t_0)\bigr)\,du
\]
gives a smooth function \(h\) with \(h(t_0)=z'(t_0)\neq0\).
Thus \(r(t)=|t-t_0|\,|h(t)|_{\mathbf g_p}\) is smooth up to the
endpoint on either side of \(t_0\). Adding these zeros to an
admissible partition gives a finite partition on which \(r\) is
piecewise smooth and continuous. Away from the partition points,
\[
 |r'(t)|\leq|\dot c(t)|_{\mathbf g_{c(t)}}.
\]
Assigning arbitrary values to the derivatives at those finitely many
points does not change their integrals. The fundamental theorem of
calculus, applied on each piece and then summed, therefore gives
\[
\begin{aligned}
 L_{\mathbf g}(c)
 &=\int_a^b|\dot c(t)|_{\mathbf g_{c(t)}}\,dt
 \geq\int_a^b|r'(t)|\,dt\\
 &\geq\left|\int_a^b r'(t)\,dt\right|
 =|r(b)-r(a)|=|v|_{\mathbf g_p}.
\end{aligned}
\]
The same conclusion is immediate for a constant curve.
By metric compatibility, \(\gamma_v\) has constant speed
\(|v|_{\mathbf g_p}\), so its length attains this bound.

Now consider an admissible curve
\(c\colon[a,b]\longrightarrow M\) that starts at \(p\) and leaves
\(U\). Its first exit time
\[
 t_*:=\min\{t\in[a,b]\mid c(t)\notin U\}
\]
satisfies \(t_*>a\). For \(a\leq t<t_*\), write
\(z(t)=\exp_p^{-1}(c(t))\).
Then \(|z(t)|_{\mathbf g_p}\longrightarrow\varepsilon\) as
\(t\to t_*^-\): otherwise, there would be \(\eta>0\) and a sequence
\(t_j\to t_*^-\) such that
\(|z(t_j)|_{\mathbf g_p}\leq\varepsilon-\eta\).
A subsequence would converge in this compact ball to \(z_*\), and
continuity would give \(c(t_*)=\exp_p(z_*)\in U\), a contradiction.
For each \(0<R<\varepsilon\), continuity yields
\(t_R<t_*\) such that \(\rho(c(t_R))=R\).
Applying the established estimate to \(c|_{[a,t_R]}\), we obtain
\[
 L_{\mathbf g}(c)\geq
 L_{\mathbf g}(c|_{[a,t_R]})\geq R.
\]
Since this holds for every \(R<\varepsilon\),
\(L_{\mathbf g}(c)\geq\varepsilon\). Thus no curve leaving \(U\)
can compete with \(\gamma_v\) when
\(|v|_{\mathbf g_p}<\varepsilon\). It also follows that every point
outside \(U\), in the component of \(p\), has distance at least
\(\varepsilon\) from \(p\). Together with the distance formula inside
\(U\), this proves \(U=B_{\mathbf g}(p,\varepsilon)\).

It remains to prove uniqueness. If \(v\neq0\) and \(c\) attains length
\(|v|_{\mathbf g_p}\), it stays in \(U\), and every inequality in the
chain of integrals is an equality. Their differences have nonnegative
integrands that are continuous on each piece. Consequently,
\[
 r'(t)=|\dot c(t)|_{\mathbf g_{c(t)}}\geq0
\]
away from the partition points. Hence \(r\) is nondecreasing.
Since \(c\) is piecewise regular, it cannot be constant on any
nondegenerate interval; therefore \(r(t)>0\) for \(t>a\).
Equality in Cauchy--Schwarz now implies
\[
 \dot c(t)=r'(t)\mathbf R(c(t)).
\]
Applying \(d(\exp_p^{-1})_{c(t)}\), we obtain
\[
 z'(t)=\frac{r'(t)}{r(t)}z(t),
 \qquad
 \frac{d}{dt}\left(\frac{z(t)}{r(t)}\right)=0
\]
on every piece contained in \((a,b]\).
Continuity across the partition points gives
\(z(t)/r(t)=v/|v|_{\mathbf g_p}\). Thus
\[
 c(t)=\gamma_v(\lambda(t)),
 \qquad
 \lambda(t):=\frac{r(t)}{|v|_{\mathbf g_p}},
 \qquad a\leq t\leq b.
\]
The function \(\lambda\) increases from \(0\) to \(1\), with positive
derivative on each piece and positive one-sided limits at the
endpoints. It is an admissible orientation-preserving
reparametrization. If \(v=0\), zero length forces the curve to be
constant, which completes the proof.
\end{proof}

\begin{lemma}[Uniformly normal and convex neighborhoods]
\label{lem:vecindades-normales-convexas-locales}
Let \((M,\mathbf g)\) be a Riemannian manifold without boundary and let
\(p\in M\). For every open neighborhood \(O\) of \(p\), there are an
open neighborhood \(W\subseteq O\) of \(p\) and \(\varepsilon>0\)
such that, for every \(x\in W\), \(\exp_x\) is a diffeomorphism on
\(B_{\mathbf g_x}(0,2\varepsilon)\), and
\[
 \exp_x\bigl(B_{\mathbf g_x}(0,2\varepsilon)\bigr)\subseteq O,
 \qquad
 W\subseteq\exp_x\bigl(B_{\mathbf g_x}(0,\varepsilon/2)\bigr).
\]
Moreover, \(p\) has a neighborhood basis of normal coordinate balls
\(V\) such that any two points of \(V\) are joined by a unique
minimizing geodesic in \(M\), parametrized on \([0,1]\), and that
geodesic stays in \(V\).
\end{lemma}
\begin{proof}
On the open set
\[
 \mathcal D_O:=\{(x,v)\in\mathcal E\mid x\in O,\ \exp_x(v)\in O\}
 \subseteq TM
\]
consider
\[
 \mathcal A\colon\mathcal D_O\longrightarrow O\times O,
 \qquad \mathcal A(x,v):=(x,\exp_x(v)).
\]
In induced coordinates on \(TM\), the identity \(\exp_x(0)=x\)
and \(d(\exp_p)_0=\operatorname{Id}\) give
\[
 d\mathcal A_{(p,0)}(h,k)=(h,h+k),
 \qquad h,k\in\mathbb R^n,\quad n=\dim M.
\]
The inverse of this linear map is \((h,\ell)\mapsto(h,\ell-h)\).
Theorem~\ref{teo: funcion inversa variedades} yields open sets
\(\mathcal U\subseteq\mathcal D_O\) and
\(\mathcal V\subseteq O\times O\), containing \((p,0)\) and
\((p,p)\), respectively, such that
\(\mathcal A\colon\mathcal U\longrightarrow\mathcal V\) is a
diffeomorphism.

Choose a smooth orthonormal frame
\((\mathbf E_1,\ldots,\mathbf E_n)\) near \(p\), as in
Proposition~\ref{marco ortonrmal existencia}. In its trivialization,
\((x,v)\) corresponds to \((x,(v^1,\ldots,v^n))\), and
\[
 v=\displaystyle\sum_{i=1}^n v^i\mathbf E_i(x),
 \qquad |v|_{\mathbf g_x}^2=\displaystyle\sum_{i=1}^n(v^i)^2.
\]
Since \(\mathcal U\) is open, we can choose an open neighborhood
\(W_0\) of \(p\) and \(\varepsilon>0\) such that
\[
 \{(x,v)\mid x\in W_0,\ |v|_{\mathbf g_x}<3\varepsilon\}
 \subseteq\mathcal U.
\]
Write \(\mathcal A^{-1}(x,y)=(x,v_{xy})\) on \(\mathcal V\).
The function \((x,y)\mapsto|v_{xy}|_{\mathbf g_x}\) is continuous
and vanishes at \((p,p)\). Thus there is an open neighborhood
\(W\subseteq W_0\) of \(p\) such that
\[
 W\times W\subseteq\mathcal V,
 \qquad |v_{xy}|_{\mathbf g_x}<\varepsilon/2
 \quad(x,y\in W).
\]
Injectivity of \(\mathcal A\) implies injectivity of each \(\exp_x\)
on the ball of radius \(2\varepsilon\). Moreover, the first component
of the differential of \(\mathcal A\) is \(d\pi\); restricting to the
vertical subspace shows that \(d(\exp_x)_v\) is invertible.
Thus \(\exp_x\) is a diffeomorphism onto its image. The other
inclusions follow from \(\mathcal U\subseteq\mathcal D_O\) and
\(\exp_x(v_{xy})=y\). By
Proposition~\ref{prop:minimizacion-radial-vecindad-normal},
\[
 \gamma_{xy}(t):=\exp_x(tv_{xy}),\qquad 0\leq t\leq1,
\]
is the unique minimizing geodesic from \(x\) to \(y\) on this
interval. The entire segment stays in \(O\), since
\(|tv_{xy}|_{\mathbf g_x}<2\varepsilon\).

To construct the convex balls, now fix a normal chart
\(z=(z^1,\ldots,z^n)\) centered at \(p\), inside any prescribed
neighborhood. Define
\[
 f:=\displaystyle\sum_{k=1}^n(z^k)^2.
\]
The covariant derivative \(\nabla df\), or Riemannian Hessian of
\(f\), has components
\[
 H_{ij}
 =\partial_i\partial_j f
  -\displaystyle\sum_{k=1}^n\Gamma^k_{ij}\partial_k f
 =2\delta_{ij}-2\displaystyle\sum_{k=1}^n\Gamma^k_{ij}z^k,
 \qquad 1\leq i,j\leq n.
\]
This matrix is symmetric and equals \(2I_n\) at \(p\).
By continuity, there is a neighborhood \(O_1\) of \(p\) within the
chart where it is positive definite. For instance, it suffices to
require \(\|H(z)-2I_n\|_{\operatorname{op}}<1\), since then
\(\langle H(z)a,a\rangle_{\mathbb R^n}\geq|a|^2\).
Apply the first part with \(O=O_1\), obtaining \(W\) and
\(\varepsilon\).

Choose \(r_0>0\) so small that
\(z^{-1}(\overline B_{\mathrm{euc}}(0,r_0))\subseteq W\).
For \(0<r\leq r_0\), let \(V_r:=\{f<r^2\}\), within this chart,
and let \(x,y\in V_r\). The geodesic \(\gamma_{xy}\) already
constructed stays in \(O_1\). Writing
\(\xi=z\circ\gamma_{xy}\), the chain rule and the geodesic equation
give
\[
\begin{aligned}
 (f\circ\gamma_{xy})''(t)
 &=\displaystyle\sum_{i,j=1}^n
 \left(\partial_i\partial_j f
 -\displaystyle\sum_{k=1}^n\Gamma^k_{ij}\partial_k f\right)
 (\gamma_{xy}(t))\,\dot\xi^i(t)\dot\xi^j(t)\\
 &\geq0.
\end{aligned}
\]
Convexity on \([0,1]\) implies
\[
 f(\gamma_{xy}(t))
 \leq(1-t)f(x)+tf(y)<r^2,\qquad 0\leq t\leq1.
\]
Hence \(\gamma_{xy}([0,1])\subseteq V_r\).
Global uniqueness has already been proved. Since \(r\) can be
arbitrarily small, these balls form a neighborhood basis at \(p\).
\end{proof}

\begin{corollary}\label{cor:curva-minimizante-es-geodesica}
On a Riemannian manifold without boundary, every admissible minimizing
curve parametrized with constant speed is a geodesic.
\end{corollary}
\begin{proof}
Every subarc of a minimizing curve is minimizing, since replacing it
by a shorter curve would reduce the total length. Fix an interior
parameter $t_0$ and choose a strongly convex normal neighborhood of
$c(t_0)$ as in Lemma~\ref{lem:vecindades-normales-convexas-locales}.
For sufficiently small $\delta>0$, the subarc
$c|_{[t_0-\delta,t_0+\delta]}$ stays in that neighborhood. Its image
must be the unique minimizing geodesic joining its endpoints.
By radial uniqueness, this subarc is $\sigma\circ\lambda$, where
$\sigma\colon[0,1]\to M$ is the geodesic joining its endpoints and
$\lambda$ preserves orientation. If its constant speed is $v>0$
and $\sigma$ has length $\ell$, then $v=\ell\lambda'$ on every
piece. Continuity and the endpoint values give
$\lambda(t)=(t-(t_0-\delta))/(2\delta)$ throughout the subinterval.
Thus the parameter is affine even across a partition point.
If the speed is zero, the curve is constant. Thus $c$ is a geodesic near every interior
parameter, and hence on its domain.
\end{proof}

\begin{definition}[Injectivity radius]\label{def: radio de inyectividad}\index{injectivity radius}
Let $(M,\mathbf{g})$ be a Riemannian manifold without boundary and let $p \in M$. The \emph{injectivity radius at $p$}, denoted by $\operatorname{inj}(p)$, is the supremum of the numbers $r>0$ for which the Riemannian exponential map
\[
\exp_p\colon B_{\mathbf{g}_p}(0,r) \subset T_pM \longrightarrow M
\]
is a diffeomorphism onto its image.

The \emph{injectivity radius of $M$} is defined by
\[
\operatorname{inj}(M) := \displaystyle\inf_{p \in M} \operatorname{inj}(p).
\]
\end{definition}

The inverse function theorem applied to $\exp_p$ at zero shows that $\operatorname{inj}(p)>0$. The supremum may be infinite, whereas the infimum over $M$ may be zero. The definition does not require completeness; the continuity result proved below will be stated under that hypothesis.

Local minimization does not guarantee minimization between distant endpoints. The following theorem relates the global existence of minimizers to completeness.

\begin{figura}[h]
    \includegraphics[width=0.7\linewidth]{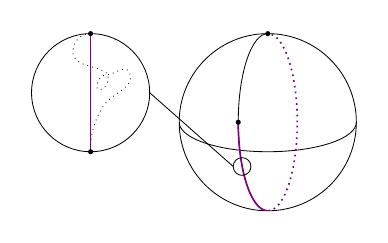}
    \caption{Every constant-speed minimizing curve is a geodesic. A geodesic is guaranteed to minimize only on sufficiently small intervals.}
    \label{fig:geodesica-minimiza-localmente}
\end{figura}
Recall that a Riemannian manifold is \emph{geodesically complete}
if every maximal geodesic is defined on $\mathbb R$.
\begin{theorem}[Hopf--Rinow]
\label{teo:variedades-riemannianas-hopf-rinow}
\index{Hopf Rinow@Hopf--Rinow}
Let $(M,\mathbf{g})$ be a connected Riemannian manifold without boundary. The following are
equivalent:
\begin{enumerate}[label=(\alph*)]
\item $(M,d_{\mathbf{g}})$ is a complete metric space;
\item every closed bounded subset of $(M,d_{\mathbf{g}})$ is compact;
\item $M$ is geodesically complete;
\item for every $p\in M$, $\exp_p$ is defined on all of $T_pM$;
\item there exists $p_0\in M$ such that $\exp_{p_0}$ is defined on all of
$T_{p_0}M$.
\end{enumerate}
If any of these conditions holds, any two points of
$M$ are joined by a geodesic segment that realizes their distance.
\end{theorem}
The proof follows Lee's strategy~\cite{LeeR}.
\begin{proof}
Conditions (c) and (d) are equivalent by the definition of
$\exp_p(v)=\gamma_v(1)$ and affine changes of parameter for
geodesics; moreover, (d) implies (e).

Assume (e), fix \(q\in M\), and set
\(r=d_{\mathbf g}(p_0,q)\). If \(r=0\), there is nothing to prove. If
\(r>0\), consider
\[
\begin{aligned}
A:=\bigl\{s\in[0,r]\ \bigm|\ &\text{there exists }v_s\in T_{p_0}M,
|v_s|_{\mathbf g}=1,\ \text{such that}\\
&x_s:=\exp_{p_0}(s v_s),\quad
d_{\mathbf g}(p_0,x_s)=s,\quad
d_{\mathbf g}(x_s,q)=r-s\bigr\}.
\end{aligned}
\]
We have \(0\in A\). The set \(A\) is closed: if \(s_j\to s\) and
\(s_j\in A\), the unit sphere of \(T_{p_0}M\) is compact, so a
subsequence of \(v_{s_j}\) converges to a unit vector \(v\).
Continuity of \(\exp_{p_0}\) and of the distance preserves the three
equalities defining \(A\).

Let \(s\in A\) with \(s<r\) and write \(x=x_s\). First choose a
normal neighborhood \(V\) of \(x\) with compact closure. This is possible by
local compactness and because normal neighborhoods form a basis at
\(x\). Since the Riemannian distance induces the manifold topology,
there exists \(\delta>0\), smaller than \(r-s\), such that
\(\overline B_{\mathbf g}(x,\delta)\subseteq V\). This ball is closed in
\(M\) and contained in the compact set \(\overline V\); therefore, it is
compact. By the definition of distance, there exist
curves \(c_j\) from \(x\) to \(q\), parametrized by arc length, with
\[
L_{\mathbf g}(c_j)\longrightarrow d_{\mathbf g}(x,q)=r-s.
\]
Let \(y_j=c_j(\delta)\). Then
\[
d_{\mathbf g}(x,y_j)\leq\delta,\qquad
d_{\mathbf g}(y_j,q)\leq L_{\mathbf g}(c_j)-\delta.
\]
By compactness, we may assume that \(y_j\to y\). Passing to the limit in the triangle inequality
\[
r-s=d_{\mathbf g}(x,q)
\leq d_{\mathbf g}(x,y_j)+d_{\mathbf g}(y_j,q)
\leq d_{\mathbf g}(x,y_j)+L_{\mathbf g}(c_j)-\delta
\]
gives
\[
d_{\mathbf g}(x,y)=\delta,\qquad
d_{\mathbf g}(y,q)=r-s-\delta.
\]
Proposition~\ref{prop:minimizacion-radial-vecindad-normal} gives
a minimizing geodesic from \(x\) to \(y\). Concatenating it with the minimizing radial
segment from \(p_0\) to \(x\) gives a curve of length
\(s+\delta\). Moreover,
\[
r\leq d_{\mathbf g}(p_0,y)+d_{\mathbf g}(y,q)
\leq s+\delta+r-s-\delta=r,
\]
so the concatenation is minimizing. If \(s=0\), choose
\(v_0\) to be the initial unit velocity of the segment from \(p_0=x\)
to \(y\); the preceding equalities give \(\delta\in A\) directly.
If \(s>0\),
Corollary~\ref{cor:curva-minimizante-es-geodesica} shows that it has no
corner at \(x\); by uniqueness for the geodesic equation, it is the extension
of \(t\mapsto\exp_{p_0}(t v_s)\). Therefore, \(s+\delta\in A\).

Since \(A\) is closed, it contains its supremum. The preceding extension
property prevents this supremum from being less than \(r\), so \(r\in A\).
Thus there exists a vector \(v\in T_{p_0}M\) with
\[
\exp_{p_0}(v)=q,\qquad
|v|_{\mathbf g}=r,
\]
and the corresponding radial geodesic is minimizing.

Now let \((q_j)\) be a Cauchy sequence. It is bounded, so
\(d_{\mathbf g}(p_0,q_j)\leq R\) for some \(R\). By the preceding argument,
we may write \(q_j=\exp_{p_0}(v_j)\) with \(|v_j|_{\mathbf g}\leq R\).
The closed ball of radius \(R\) in the finite-dimensional space
\(T_{p_0}M\) is compact; a subsequence \(v_{j_\ell}\) converges to \(v\),
and then \(q_{j_\ell}\to\exp_{p_0}(v)\). Since the original sequence is
Cauchy, the entire sequence converges to the same point. This proves
\((e)\Rightarrow(a)\).

Now assume (a). If a maximal geodesic
$\gamma\colon[0,b)\to M$ had $b<\infty$, its speed would be constant and
\[
 d_{\mathbf{g}}(\gamma(s),\gamma(t))
 \leq L_{\mathbf{g}}(\gamma\restriction_{[s,t]})
 =|t-s|\,|\dot\gamma(0)|_{\mathbf{g}}.
\]
Thus, $\gamma(t)$ would have a limit $p\in M$ as $t\to b^{-}$.
Choose a chart around $p$ and a smaller neighborhood whose compact
closure is contained in that chart. A tail of $\gamma$ stays in
the smaller neighborhood. On its closure, the metric is uniformly equivalent
to the Euclidean metric and the Christoffel symbols are bounded.
If $x^i(t)$ are the coordinates of $\gamma(t)$ and $v^i(t)=\dot x^i(t)$,
the constant Riemannian speed bounds all the components $v^i(t)$.
The geodesic equation
\[
 \dot v^k(t)=-\Gamma^k_{ij}(x(t))v^i(t)v^j(t)
\]
then bounds their derivatives. Hence, $v(t)$ is uniformly
Lipschitz on this tail and has a limit $v_b$ as $t\to b^{-}$.
Extend $(x(t),v(t))$ continuously to the data $(x(p),v_b)$.
The integral equation for the geodesic system passes to the limit at $b$;
local existence and uniqueness for this system, with these initial
data, extend the solution beyond $b$, a contradiction. The same
argument at the left endpoint proves (c). We now also have (d);
applying the radial argument of the preceding paragraph with any point
as center gives a minimizing geodesic segment between two arbitrary
points.

Under these conditions, for $p\in M$ and $R>0$, the existence of
minimizers gives
\[
 \overline B_{\mathbf{g}}(p,R)
 =\exp_p\bigl(\overline B_{\mathbf{g}_p}(0,R)\bigr).
\]
The right-hand side is the continuous image of a compact
finite-dimensional set, so every closed ball is compact. Every closed bounded
set is contained as a closed subset in one of these balls, giving
(b). Finally, (b) implies (a): a Cauchy sequence is bounded,
has a convergent subsequence in a compact closed ball, and,
being Cauchy, converges to the same limit. This proves all the
equivalences and consequences.
\end{proof}

We record separately three consequences that will be used later:
\begin{corollary}\label{cor:variedades-riemannianas-variedad-riemanniana-conexa-punto-tal-mapeo}
If $(M,\mathbf{g})$ is a connected Riemannian manifold without boundary and there exists
$p\in M$ such that $\exp_p$ is defined on all of $T_pM$, then $M$ is
complete.
\end{corollary}
\begin{proof}
This is implication \((e)\Rightarrow(a)\) of
Theorem~\ref{teo:variedades-riemannianas-hopf-rinow}, taking \(p_0=p\).
\end{proof}
\begin{corollary}\label{dos puntos se pueden unir por segmento geodésico minimizante}
 If $M$ is a connected complete Riemannian manifold without boundary,
 then any two points of $M$ can be joined by a minimizing
 geodesic segment.
\end{corollary}
\begin{proof}
Metric completeness is condition (a) of
Hopf--Rinow; the last assertion of the theorem gives the minimizing
segment.
\end{proof}
\begin{corollary}\label{cor:variedades-riemannianas-variedad-riemanniana-conexa-geodesica-maximal-esta}
 If $M$ is a connected complete Riemannian manifold without boundary,
 then every maximal geodesic in $M$ is defined on all of $\mathbb R$.
\end{corollary}
\begin{proof}
By Hopf--Rinow, metric completeness is equivalent to geodesic
completeness, which is precisely the assertion.
\end{proof}
\section{Curvature}
Curvature is the central notion in the study of Riemannian geometry. Intuitively, it measures the failure of covariant derivatives to commute: transporting a vector in different directions and then comparing the resulting vectors reveals the presence of curvature.

\begin{definition}\label{def:variedades-riemannianas-variedad-5}\index{Riemann curvature tensor}
 Let $(M,\mathbf{g})$ be a Riemannian manifold with or without boundary. The $(1,3)$-\textbf{\textit{Riemann curvature tensor}} of $M$ is the tensor field $\mathbf{R}\in \Gamma(T^{(1,3)}(TM))$ of type $(1,3)$ given by

 \[\mathbf{R}(\mathbf{X},\mathbf{Y})(\mathbf{Z}):=
 \nabla_{\mathbf{X}}\nabla _{\mathbf{Y}}\mathbf{Z}
 -\nabla _{\mathbf{Y}}\nabla_{\mathbf{X}}\mathbf{Z}
 -\nabla_{[\mathbf{X},\mathbf{Y}]}\mathbf{Z}.\]
\end{definition}

\begin{proposition}\label{prop:variedades-riemannianas-tensor-curvatura-riemann-define-campo-tensorial}
Let $(M,\mathbf{g})$ be a Riemannian manifold with or without boundary. The
$(1,3)$-Riemann curvature tensor $\mathbf{R}$ defines a tensor field
of type $(1,3)$ on $M$.
\end{proposition}

\begin{proof}
By Lemma~\ref{lema de caracterización tensorial}, it suffices to verify that
$\mathbf{R}$ is $C^{\infty}(M)$-multilinear in each argument. Additivity
in all three arguments follows from the linearity of the connection and the
bracket; it remains to check compatibility with multiplication by a function.
Take $f\in C^{\infty}(M)$ and
$\mathbf{X},\mathbf{Y},\mathbf{Z}\in\mathfrak X(M)$.

Let us verify this for the first argument:
\[
\begin{aligned}
\mathbf{R}(f\mathbf{X},\mathbf{Y})\mathbf{Z}
&= \nabla_{f\mathbf{X}}\nabla_{\mathbf{Y}}\mathbf{Z}
- \nabla_{\mathbf{Y}}\nabla_{f\mathbf{X}}\mathbf{Z}
- \nabla_{[f\mathbf{X},\mathbf{Y}]}\mathbf{Z} \\
&= f\nabla_{\mathbf{X}}\nabla_{\mathbf{Y}}\mathbf{Z}
- \nabla_{\mathbf{Y}}(f\nabla_{\mathbf{X}}\mathbf{Z})
- \nabla_{f[\mathbf{X},\mathbf{Y}]-(\mathbf{Y}f)\mathbf{X}}\mathbf{Z} \\
&= f\nabla_{\mathbf{X}}\nabla_{\mathbf{Y}}\mathbf{Z}
-(\mathbf{Y}f)\nabla_{\mathbf{X}}\mathbf{Z}
-f\nabla_{\mathbf{Y}}\nabla_{\mathbf{X}}\mathbf{Z}
-f\nabla_{[\mathbf{X},\mathbf{Y}]}\mathbf{Z}
+(\mathbf{Y}f)\nabla_{\mathbf{X}}\mathbf{Z} \\
&= f\mathbf{R}(\mathbf{X},\mathbf{Y})\mathbf{Z}.
\end{aligned}
\]

For the second argument:
\[
\begin{aligned}
\mathbf{R}(\mathbf{X},f\mathbf{Y})\mathbf{Z}
&= \nabla_{\mathbf{X}}\nabla_{f\mathbf{Y}}\mathbf{Z}
- \nabla_{f\mathbf{Y}}\nabla_{\mathbf{X}}\mathbf{Z}
- \nabla_{[\mathbf{X},f\mathbf{Y}]}\mathbf{Z} \\
&= \nabla_{\mathbf{X}}(f\nabla_{\mathbf{Y}}\mathbf{Z})
- f\nabla_{\mathbf{Y}}\nabla_{\mathbf{X}}\mathbf{Z}
- \nabla_{f[\mathbf{X},\mathbf{Y}]+(\mathbf{X}f)\mathbf{Y}}\mathbf{Z} \\
&= (\mathbf{X}f)\nabla_{\mathbf{Y}}\mathbf{Z}
+f\nabla_{\mathbf{X}}\nabla_{\mathbf{Y}}\mathbf{Z}
-f\nabla_{\mathbf{Y}}\nabla_{\mathbf{X}}\mathbf{Z}
-f\nabla_{[\mathbf{X},\mathbf{Y}]}\mathbf{Z}
-(\mathbf{X}f)\nabla_{\mathbf{Y}}\mathbf{Z} \\
&= f\mathbf{R}(\mathbf{X},\mathbf{Y})\mathbf{Z}.
\end{aligned}
\]

Finally, for the third and last argument:
\[
\begin{aligned}
\mathbf{R}(\mathbf{X},\mathbf{Y})(f\mathbf{Z})
&= \nabla_{\mathbf{X}}\nabla_{\mathbf{Y}}(f\mathbf{Z})
- \nabla_{\mathbf{Y}}\nabla_{\mathbf{X}}(f\mathbf{Z})
- \nabla_{[\mathbf{X},\mathbf{Y}]}(f\mathbf{Z}) \\
&= \nabla_{\mathbf{X}}\big((\mathbf{Y}f)\mathbf{Z}
+ f\nabla_{\mathbf{Y}}\mathbf{Z}\big)
- \nabla_{\mathbf{Y}}\big((\mathbf{X}f)\mathbf{Z}
+ f\nabla_{\mathbf{X}}\mathbf{Z}\big)
- ([\mathbf{X},\mathbf{Y}]f)\mathbf{Z}
- f\nabla_{[\mathbf{X},\mathbf{Y}]}\mathbf{Z} \\
&= (\mathbf{X}(\mathbf{Y}f))\mathbf{Z}
+ (\mathbf{Y}f)\nabla_{\mathbf{X}}\mathbf{Z}
+ (\mathbf{X}f)\nabla_{\mathbf{Y}}\mathbf{Z}
+ f\nabla_{\mathbf{X}}\nabla_{\mathbf{Y}}\mathbf{Z} \\
&\quad - (\mathbf{Y}(\mathbf{X}f))\mathbf{Z}
- (\mathbf{X}f)\nabla_{\mathbf{Y}}\mathbf{Z}
- (\mathbf{Y}f)\nabla_{\mathbf{X}}\mathbf{Z}
- f\nabla_{\mathbf{Y}}\nabla_{\mathbf{X}}\mathbf{Z} \\
&\quad - ([\mathbf{X},\mathbf{Y}]f)\mathbf{Z}
- f\nabla_{[\mathbf{X},\mathbf{Y}]}\mathbf{Z}.
\end{aligned}
\]
The terms containing first derivatives of $f$ cancel. Moreover,
\[
\mathbf{X}(\mathbf{Y}f) - \mathbf{Y}(\mathbf{X}f) - [\mathbf{X},\mathbf{Y}]f = 0,
\]
by the definition of the Lie bracket acting on functions. Consequently,
\[
\mathbf{R}(\mathbf{X},\mathbf{Y})(f\mathbf{Z})
= f\mathbf{R}(\mathbf{X},\mathbf{Y})\mathbf{Z}.
\]

We have verified multilinearity in all three arguments. Therefore,
$\mathbf{R}$ is a tensor field of type $(1,3)$ on $M$.
\end{proof}
Since $\mathbf{R}$ is a tensor field of type $(1,3)$, we may write it
in local coordinates $(x^{i})$ as
\[
\mathbf{R}=R_{ijk}{}^{l}\,\mathbf{d}x^{i}\otimes \mathbf{d}x^{j}\otimes \mathbf{d}x^{k}\otimes
\boldsymbol{\partial}_{l}.
\]
The small space before the index $l$ visually separates the three
covariant indices from the contravariant index.

The coefficients satisfy the relation
\[
\mathbf{R}(\boldsymbol{\partial}_{i},\boldsymbol{\partial}_{j})
\boldsymbol{\partial}_{k}=R_{ijk}{}^{l}\boldsymbol{\partial}_{l}.
\]

The following proposition shows how to compute the components of $R$ in coordinates:
\begin{proposition}\label{1,3 curvatura en coordenadas}
 Let $(M,\mathbf{g})$ be a Riemannian manifold with or without boundary. The components of the $(1,3)$-curvature tensor in coordinates $(x^{i})$ are given by \[R_{ijk}{}^{l}=\partial_i\Gamma_{jk}^{l}-\partial_{j}\Gamma_{ik}^{l}+\Gamma_{jk}^{m}\Gamma_{im}^{l}-\Gamma_{ik}^{m}\Gamma_{jm}^{l}.\]
\end{proposition}
\begin{proof}
Let $(x^{i})$ be local coordinates on $M$ and let $\{\boldsymbol{\partial}_{i}\}$ be the associated coordinate frame. By the definition of the Christoffel symbols,
\[
\nabla_{\boldsymbol{\partial}_{i}}\boldsymbol{\partial}_{j}=\Gamma^{m}_{ij}\boldsymbol{\partial}_{m}.
\]
Using the linearity and Leibniz rule of the connection, we have
\[
\nabla_{\boldsymbol{\partial}_{i}}\nabla_{\boldsymbol{\partial}_{j}}\boldsymbol{\partial}_{k}
=\nabla_{\boldsymbol{\partial}_{i}}\big(\Gamma^{l}_{jk}\boldsymbol{\partial}_{l}\big)
=\partial_{i}\Gamma^{l}_{jk}\boldsymbol{\partial}_{l}
+\Gamma^{l}_{jk}\nabla_{\boldsymbol{\partial}_{i}}\boldsymbol{\partial}_{l}
=\partial_{i}\Gamma^{l}_{jk}\boldsymbol{\partial}_{l}
+\Gamma^{l}_{jk}\Gamma^{m}_{il}\boldsymbol{\partial}_{m}.
\]
Interchanging the indices $i$ and $j$ in the preceding calculation gives
\[
\nabla_{\boldsymbol{\partial}_{j}}\nabla_{\boldsymbol{\partial}_{i}}\boldsymbol{\partial}_{k}
=\partial_{j}\Gamma^{l}_{ik}\boldsymbol{\partial}_{l}
+\Gamma^{l}_{ik}\Gamma^{m}_{jl}\boldsymbol{\partial}_{m}.
\]
Since $[\boldsymbol{\partial}_{i},\boldsymbol{\partial}_{j}]=0$ in our coordinate frame, the third term in the definition vanishes:
\[
\nabla_{[\boldsymbol{\partial}_{i},\boldsymbol{\partial}_{j}]}\boldsymbol{\partial}_{k}=0.
\]
Therefore, recalling that
$\mathbf{R}(\boldsymbol{\partial}_{i},\boldsymbol{\partial}_{j})
\boldsymbol{\partial}_{k}=\nabla_{\boldsymbol{\partial}_{i}}
\nabla_{\boldsymbol{\partial}_{j}}\boldsymbol{\partial}_{k}
-\nabla_{\boldsymbol{\partial}_{j}}\nabla_{\boldsymbol{\partial}_{i}}
\boldsymbol{\partial}_{k}
-\nabla_{[\boldsymbol{\partial}_{i},\boldsymbol{\partial}_{j}]}
\boldsymbol{\partial}_{k}$ and interchanging $l$ and $m$ in
$\Gamma_{jk}^{l}\Gamma_{il}^{m}\boldsymbol{\partial}_{m}
-\Gamma_{ik}^{l}\Gamma_{jl}^{m}\boldsymbol{\partial}_{m}$,
\[
\mathbf{R}(\boldsymbol{\partial}_{i},\boldsymbol{\partial}_{j})\boldsymbol{\partial}_{k}
=\big(\partial_{i}\Gamma^{l}_{jk}-\partial_{j}\Gamma^{l}_{ik}
+\Gamma^{m}_{jk}\Gamma^{l}_{im}
-\Gamma^{m}_{ik}\Gamma^{l}_{jm}\big)\boldsymbol{\partial}_{l}.
\]
Comparing coefficients in the basis $\{\boldsymbol{\partial}_{l}\}$, we conclude that
\[
R_{ijk}{}^{l}
=\partial_{i}\Gamma^{l}_{jk}-\partial_{j}\Gamma^{l}_{ik}
+\Gamma^{m}_{jk}\Gamma^{l}_{im}
-\Gamma^{m}_{ik}\Gamma^{l}_{jm},
\]
which is the desired formula.
\end{proof}

For any vector fields
$\mathbf{X},\mathbf{Y}\in \mathfrak{X}(M)$, the map
$\mathbf{R}(\mathbf{X},\mathbf{Y})\colon \mathfrak{X}(M)\longrightarrow
\mathfrak{X}(M)$ given by
$\mathbf{Z}\mapsto \mathbf{R}(\mathbf{X},\mathbf{Y})\mathbf{Z}$ is a
smooth homomorphism of the tangent bundle; in particular,
$\mathbf{R}(\mathbf{X},\mathbf{Y})$ is $C^\infty(M)$-linear. It is called the
\textit{\textbf{curvature endomorphism determined by $\mathbf{X}$ and
$\mathbf{Y}$}}. Lowering the contravariant index with the metric gives
the covariant tensor field defined next.
\begin{definition}\label{def:variedades-riemannianas-tensor-de-curvatura-de-riemann}\index{Riemann curvature tensor}
 Let $(M,\mathbf{g})$ be a Riemannian manifold with or without boundary. Define
 the \textit{\textbf{Riemann curvature tensor}} to be the tensor field
 of type $(0,4)$ given by
 \[
 \mathbf{Rm}(\mathbf{X},\mathbf{Y},\mathbf{Z},\mathbf{W})
 =\langle \mathbf{R}(\mathbf{X},\mathbf{Y})\mathbf{Z},\mathbf{W}\rangle_{\mathbf{g}}.
 \]
\end{definition}
As a corollary of Proposition~\ref{1,3 curvatura en coordenadas}, we obtain a local expression for the Riemann curvature tensor:
\begin{corollary}\label{Rm en coordenadas}
 If $(M,\mathbf{g})$ is a Riemannian manifold with or without boundary and
 \[
 \mathbf{Rm}=R_{ijkl}\,\mathbf{d}x^{i}\otimes \mathbf{d}x^{j}\otimes \mathbf{d}x^{k}\otimes \mathbf{d}x^{l},
 \]
 then
 \[
 R_{ijkl}=g_{lm}\left(\partial_{i}\Gamma_{jk}^{m}
 -\partial_{j}\Gamma_{ik}^{m}+\Gamma_{jk}^{p}\Gamma_{ip}^{m}
 -\Gamma_{ik}^{p}\Gamma_{jp}^{m}\right).
 \]
\end{corollary}
\begin{proof}
By the definition of $\mathbf{Rm}$ and the formula in
Proposition~\ref{1,3 curvatura en coordenadas},
\[
R_{ijkl}
=
\mathbf{Rm}(\boldsymbol{\partial}_i,\boldsymbol{\partial}_j,
\boldsymbol{\partial}_k,\boldsymbol{\partial}_l)
=
g_{lm}R^{m}{}_{ijk}.
\]
Substituting
\[
R^{m}{}_{ijk}
=
\partial_i\Gamma^m_{jk}-\partial_j\Gamma^m_{ik}
+\Gamma^p_{jk}\Gamma^m_{ip}
-\Gamma^p_{ik}\Gamma^m_{jp}
\]
gives the stated expression.
\end{proof}
The symmetries and antisymmetries of the Riemann curvature tensor relate its components and simplify calculations involving it.
\begin{proposition}[Symmetries of the curvature tensor]\label{simetrías intrínsecas del tensor de curvatura}\index{symmetries of the curvature tensor}
Let $(M,\mathbf{g})$ be a Riemannian manifold with or without boundary. The Riemann curvature tensor has the following properties for all $\mathbf{X},\mathbf{Y},\mathbf{Z},\mathbf{W}\in \mathfrak{X}(M)$.
\begin{enumerate}
 \item $\mathbf{Rm}(\mathbf{W},\mathbf{X},\mathbf{Y},\mathbf{Z})=-\mathbf{Rm}(\mathbf{X},\mathbf{W},\mathbf{Y},\mathbf{Z})$.
 \item $\mathbf{Rm}(\mathbf{W},\mathbf{X},\mathbf{Y},\mathbf{Z})=-\mathbf{Rm}(\mathbf{W},\mathbf{X},\mathbf{Z},\mathbf{Y})$.
 \item $\mathbf{Rm}(\mathbf{W},\mathbf{X},\mathbf{Y},\mathbf{Z})=\mathbf{Rm}(\mathbf{Y},\mathbf{Z},\mathbf{W},\mathbf{X})$.
 \item $\mathbf{Rm}(\mathbf{W},\mathbf{X},\mathbf{Y},\mathbf{Z})+\mathbf{Rm}(\mathbf{X},\mathbf{Y},\mathbf{W},\mathbf{Z})+\mathbf{Rm}(\mathbf{Y},\mathbf{W},\mathbf{X},\mathbf{Z})=0$. This identity is known as the \textit{\textbf{Bianchi identity}}.
\end{enumerate}
\end{proposition}
\begin{proof}
Write
\[
R(\mathbf X,\mathbf Y,\mathbf Z,\mathbf W)
:=\langle\mathbf R(\mathbf X,\mathbf Y)\mathbf Z,\mathbf W\rangle.
\]
The definition of \(\mathbf R\) immediately gives
\(\mathbf R(\mathbf X,\mathbf Y)=-\mathbf R(\mathbf Y,\mathbf X)\),
proving (1).

Metric compatibility of the connection implies
\[
\begin{aligned}
&R(\mathbf X,\mathbf Y,\mathbf Z,\mathbf W)
+R(\mathbf X,\mathbf Y,\mathbf W,\mathbf Z)\\
&\quad=
\mathbf X\mathbf Y\langle\mathbf Z,\mathbf W\rangle
-\mathbf Y\mathbf X\langle\mathbf Z,\mathbf W\rangle
-[\mathbf X,\mathbf Y]\langle\mathbf Z,\mathbf W\rangle=0.
\end{aligned}
\]
In the first equality, the three covariant derivatives in each
curvature term are expanded and grouped using the compatibility rule; the last
equality is the definition of the bracket acting on functions. This
proves (2).

Since the connection is torsion-free, expanding the cyclic sum gives
\[
\begin{aligned}
&\mathbf R(\mathbf X,\mathbf Y)\mathbf Z
+\mathbf R(\mathbf Y,\mathbf Z)\mathbf X
+\mathbf R(\mathbf Z,\mathbf X)\mathbf Y\\
&=\nabla_{\mathbf X}
 (\nabla_{\mathbf Y}\mathbf Z-\nabla_{\mathbf Z}\mathbf Y)
+\nabla_{\mathbf Y}
 (\nabla_{\mathbf Z}\mathbf X-\nabla_{\mathbf X}\mathbf Z)
+\nabla_{\mathbf Z}
 (\nabla_{\mathbf X}\mathbf Y-\nabla_{\mathbf Y}\mathbf X)\\
&\quad-\nabla_{[\mathbf X,\mathbf Y]}\mathbf Z
-\nabla_{[\mathbf Y,\mathbf Z]}\mathbf X
-\nabla_{[\mathbf Z,\mathbf X]}\mathbf Y\\
&=\nabla_{\mathbf X}[\mathbf Y,\mathbf Z]
+\nabla_{\mathbf Y}[\mathbf Z,\mathbf X]
+\nabla_{\mathbf Z}[\mathbf X,\mathbf Y]\\
&\quad-\nabla_{[\mathbf X,\mathbf Y]}\mathbf Z
-\nabla_{[\mathbf Y,\mathbf Z]}\mathbf X
-\nabla_{[\mathbf Z,\mathbf X]}\mathbf Y=0.
\end{aligned}
\]
In the last equality, terms are grouped, for example, as
\(\nabla_{\mathbf X}[\mathbf Y,\mathbf Z]
-\nabla_{[\mathbf Y,\mathbf Z]}\mathbf X
=[\mathbf X,[\mathbf Y,\mathbf Z]]\), and the three resulting expressions
sum to zero by the Jacobi identity. Pairing with \(\mathbf W\)
gives (4).

It remains to prove (3). Fix four vector fields and abbreviate their names by
\(a,b,c,d\). Using only the antisymmetries already proved, the
Bianchi identities corresponding to the relevant permutations
can be written as
\[
\begin{aligned}
R_{bcad}&=R_{acbd}-R_{abcd},\\
R_{bdac}&=R_{abcd}+R_{adbc},\\
R_{cdab}&=R_{acbd}-R_{adbc},\\
R_{bcad}&=R_{bdac}-R_{cdab},
\end{aligned}
\]
where \(R_{abcd}:=R(a,b,c,d)\). Substituting the first three equalities
into the fourth gives
\[
R_{abcd}=R_{acbd}-R_{adbc}=R_{cdab}.
\]
This is the pair-interchange symmetry in (3).
\end{proof}
We may express the symmetries of the curvature tensor in coordinates as follows:
\begin{corollary}\label{simetrías tensor curvatura coordenadas}
 Let $(M,\mathbf{g})$ be a Riemannian manifold with or without boundary. In
 any local coordinates, the components of $\mathbf{Rm}$ satisfy:
 \begin{enumerate}
 \item $R_{ijkl}=-R_{jikl}$.
 \item $R_{ijkl}=-R_{ijlk}$.
 \item $R_{ijkl}=R_{klij}$.
 \item $R_{ijkl}+R_{jkil}+R_{kijl}=0$.
 \end{enumerate}
\end{corollary}
\begin{proof}
Apply the four parts of the preceding proposition to the coordinate
vector fields
\(\boldsymbol\partial_i,\boldsymbol\partial_j,
\boldsymbol\partial_k,\boldsymbol\partial_l\).
\end{proof}

\begin{definition}[Ricci tensor and Ricci endomorphism]
\label{def:tensor-endomorfismo-ricci}
\index{Ricci tensor}
\index{Ricci endomorphism}
Let $(M,\mathbf{g})$ be a Riemannian manifold of dimension $n$ with or without boundary. For
$X,Y\in T_pM$, define
\begin{equation}
\label{eq:definicion-tensor-ricci}
\mathbf{Ric}_p(X,Y)
:=
\displaystyle\sum_{i=1}^{n}
 \mathbf{g}_p\bigl(\mathbf{R}(\mathbf{e}_i,X)Y,\mathbf{e}_i\bigr),
\end{equation}
where $(\mathbf{e}_1,\dots,\mathbf{e}_n)$ is an orthonormal basis of $T_pM$. This expression
is independent of the chosen orthonormal basis and defines a symmetric tensor
field of type $(0,2)$, called the \emph{Ricci tensor}.

The \emph{Ricci endomorphism} is the field
$\mathbf{Ric}^{\sharp}\in\Gamma(\operatorname{End}(TM))$
determined by
\begin{equation}
\label{eq:definicion-endomorfismo-ricci}
\mathbf{g}(\mathbf{Ric}^{\sharp}X,Y)=\mathbf{Ric}(X,Y).
\end{equation}
\end{definition}

\begin{proposition}[Curvature contraction and the Ricci endomorphism]
\label{prop:contraccion-curvatura-endomorfismo-ricci}
Let $(M,\mathbf{g})$ be a Riemannian manifold with or without boundary of
dimension $n$. For every $p\in M$, every orthonormal basis
$(\mathbf{e}_1,\dots,\mathbf{e}_n)$ of $T_pM$, and every
$X\in T_pM$,
\begin{equation}
\label{eq:contraccion-curvatura-ricci-cap2}
-\displaystyle\sum_{i=1}^{n}\mathbf{R}(\mathbf{e}_i,X)\mathbf{e}_i
=
\mathbf{Ric}^{\sharp}X.
\end{equation}
In particular, $\mathbf{Ric}^{\sharp}$ is self-adjoint with respect to
$\mathbf{g}$.
\end{proposition}

\begin{proof}
Independence of the basis in
\eqref{eq:definicion-tensor-ricci} follows from basis independence of the trace of an
endomorphism. Let us verify symmetry. If $X,Y\in T_pM$, the
pair symmetries and antisymmetries of the curvature tensor give
\[
\begin{split}
\mathbf{Ric}(Y,X)
&=
\displaystyle\sum_{i=1}^{n}\mathbf{Rm}(\mathbf{e}_i,Y,X,\mathbf{e}_i)\\
&=
\displaystyle\sum_{i=1}^{n}\mathbf{Rm}(X,\mathbf{e}_i,\mathbf{e}_i,Y)\\
&=
-\displaystyle\sum_{i=1}^{n}\mathbf{Rm}(\mathbf{e}_i,X,\mathbf{e}_i,Y)\\
&=
\displaystyle\sum_{i=1}^{n}\mathbf{Rm}(\mathbf{e}_i,X,Y,\mathbf{e}_i)
=
\mathbf{Ric}(X,Y).
\end{split}
\]
To prove
\eqref{eq:contraccion-curvatura-ricci-cap2}, take $Y\in T_pM$ and use
the antisymmetry of $\mathbf{Rm}$ in its last two arguments:
\[
\begin{split}
\mathbf{g}\left(-\displaystyle\sum_{i=1}^{n}\mathbf{R}(\mathbf{e}_i,X)\mathbf{e}_i,Y\right)
&=
-\displaystyle\sum_{i=1}^{n}\mathbf{Rm}(\mathbf{e}_i,X,\mathbf{e}_i,Y)\\
&=
\displaystyle\sum_{i=1}^{n}\mathbf{Rm}(\mathbf{e}_i,X,Y,\mathbf{e}_i)
=
\mathbf{Ric}(X,Y).
\end{split}
\]
Nondegeneracy of $\mathbf{g}$ gives the identity. Finally, symmetry of
$\mathbf{Ric}$ implies
\[
\mathbf{g}(\mathbf{Ric}^{\sharp}X,Y)
=
\mathbf{Ric}(X,Y)
=
\mathbf{Ric}(Y,X)
=
\mathbf{g}(X,\mathbf{Ric}^{\sharp}Y),
\]
so $\mathbf{Ric}^{\sharp}$ is self-adjoint.
\end{proof}

\begin{definition}[Scalar curvature]
\label{def:curvatura-escalar-riemanniana}
\index{scalar curvature}
Let $(M,\mathbf{g})$ be a Riemannian manifold with or without boundary. The \emph{scalar curvature} is the
smooth function
\[
R_{\mathbf{g}}:=\operatorname{tr}_{\mathbf{g}}(\mathbf{Ric}_{\mathbf{g}}).
\]
If $(\mathbf{e}_1,\ldots,\mathbf{e}_n)$ is a local orthonormal frame, then
\[
R_{\mathbf{g}}=\displaystyle\sum_{i=1}^n\mathbf{Ric}_{\mathbf{g}}(\mathbf{e}_i,\mathbf{e}_i),
\]
and in coordinates, $R_{\mathbf{g}}=g^{ij}\operatorname{Ric}_{ij}$. In the analytic
chapters, we retain the notation $R_{\mathbf{g}}$ to distinguish the scalar from
the curvature tensor or endomorphism.
\end{definition}

\begin{definition}[Sectional curvature]\label{def:variedades-riemannianas-curvatura-seccional}\index{sectional curvature}
Let $(M,\mathbf{g})$ be a Riemannian manifold with or without boundary and let $p \in M$. Given a two-dimensional plane $\Pi \subseteq T_pM$ spanned by two linearly independent vectors $u, v \in T_pM$, the \emph{sectional curvature} of $\Pi$ is defined by
\[
K(\Pi) = \displaystyle\frac{\mathbf{g}(\mathbf{R}(u,v)v, u)}{\mathbf{g}(u,u)\mathbf{g}(v,v) - \mathbf{g}(u,v)^2},
\]
where $\mathbf{R}$ is the Riemann curvature tensor. This value depends only on the plane $\Pi$ and is independent of the basis $\{u,v\}$ chosen to span it.
\end{definition}

\begin{theorem}[Curvature tensor for spaces of constant curvature]
\label{thm: tensor curvatura constante}\index{curvature tensor for spaces of constant curvature}
Let $(M,\mathbf{g})$ be a Riemannian manifold with or without boundary. If
$(M,\mathbf{g})$ has constant sectional curvature $K$, then, for
every $p\in M$ and all $u,v,w\in T_pM$, we have
\[
\mathbf{R}(u,v)w = K \big( \mathbf{g}(v,w)u - \mathbf{g}(u,w)v \big).
\]
\end{theorem}

\begin{proof}
By hypothesis, the sectional curvature is a constant $K$. From the definition of sectional curvature, we rewrite the equation to obtain
\[
\mathbf{g}(\mathbf{R}(u,v)v, u) = K \big( \mathbf{g}(u,u)\mathbf{g}(v,v) - \mathbf{g}(u,v)^2 \big).
\]
This equality holds for every pair of linearly independent vectors $u,v$ and extends trivially to the case in which they are linearly dependent (since both sides vanish).

Given arbitrary $u,v,w$, we compute the expression
$\mathbf{g}(\mathbf{R}(u+w,v)v,u+w)$ in two different ways.

First, applying the preceding formula directly to the pair $(u+w, v)$, we obtain:
\begin{align*}
\mathbf{g}(\mathbf{R}(u+w, v)v, u+w) &= K \big( \mathbf{g}(v,v)\mathbf{g}(u+w, u+w) - \mathbf{g}(u+w, v)^2 \big) \\
&= K \big( \mathbf{g}(v,v)[\mathbf{g}(u,u) + \mathbf{g}(w,w) + 2\mathbf{g}(u,w)] \\
&\qquad - [\mathbf{g}(u,v)^2 + \mathbf{g}(w,v)^2 + 2\mathbf{g}(u,v)\mathbf{g}(w,v)] \big).
\end{align*}

Second, using multilinearity of the Riemann tensor on the left-hand side:
\begin{align*}
\mathbf{g}(\mathbf{R}(u+w, v)v, u+w) &= \mathbf{g}(\mathbf{R}(u,v)v + \mathbf{R}(w,v)v, u+w) \\
&= \mathbf{g}(\mathbf{R}(u,v)v, u) + \mathbf{g}(\mathbf{R}(w,v)v, w) + \mathbf{g}(\mathbf{R}(u,v)v, w) + \mathbf{g}(\mathbf{R}(w,v)v, u).
\end{align*}
By Proposition~\ref{simetrías intrínsecas del tensor de curvatura}
(parts 1, 2, and 3), we know that
$\mathbf{Rm}(u,v,v,w)=\mathbf{Rm}(w,v,v,u)$, which is equivalent to
$\mathbf{g}(\mathbf{R}(u,v)v,w)=\mathbf{g}(\mathbf{R}(w,v)v,u)$. Moreover,
applying our basic formula to the first two terms simplifies the expression
to
\[
K \big( \mathbf{g}(u,u)\mathbf{g}(v,v) - \mathbf{g}(u,v)^2 \big) + K \big( \mathbf{g}(w,w)\mathbf{g}(v,v) - \mathbf{g}(w,v)^2 \big) + 2\mathbf{g}(\mathbf{R}(u,v)v, w).
\]

Equating the two expansions of $\mathbf{g}(\mathbf{R}(u+w,v)v,u+w)$ and
canceling the common terms gives
\[
2\mathbf{g}(\mathbf{R}(u,v)v, w) = 2K \big( \mathbf{g}(v,v)\mathbf{g}(u,w) - \mathbf{g}(u,v)\mathbf{g}(w,v) \big).
\]
After division by $2$, this equality can be written as
\[
\mathbf{g}\bigl(\mathbf{R}(u,v)v
-K(\mathbf{g}(v,v)u-\mathbf{g}(u,v)v),w\bigr)=0
\qquad\text{for every }w\in T_pM.
\]
Nondegeneracy of $\mathbf{g}_p$ forces the first argument to be
zero. Therefore,
\begin{equation}
\label{eq: curvatura constante 1}
\mathbf{R}(u,v)v = K \big( \mathbf{g}(v,v)u - \mathbf{g}(u,v)v \big).
\end{equation}

We now polarize in the second and third arguments. Take arbitrary $u,v,w$
and compute $\mathbf{R}(u,v+w)(v+w)$ in two ways. First,
using \eqref{eq: curvatura constante 1}:
\begin{align*}
\mathbf{R}(u, v+w)(v+w) &= K \big( \mathbf{g}(v+w,v+w)u - \mathbf{g}(u,v+w)(v+w) \big) \\
&= K \big( [\mathbf{g}(v,v)+\mathbf{g}(w,w)+2\mathbf{g}(v,w)]u \\
&\qquad - \mathbf{g}(u,v)v - \mathbf{g}(u,w)v - \mathbf{g}(u,v)w - \mathbf{g}(u,w)w \big).
\end{align*}

Second, by multilinearity of $\mathbf{R}$:
\[
\mathbf{R}(u, v+w)(v+w) = \mathbf{R}(u,v)v + \mathbf{R}(u,w)w + \mathbf{R}(u,v)w + \mathbf{R}(u,w)v.
\]
Applying (\ref{eq: curvatura constante 1}) to the first two terms of this sum:
\[
K \big( \mathbf{g}(v,v)u - \mathbf{g}(u,v)v \big) + K \big( \mathbf{g}(w,w)u - \mathbf{g}(u,w)w \big) + \mathbf{R}(u,v)w + \mathbf{R}(u,w)v.
\]

Equating these two new expansions and canceling corresponding terms, we arrive at:
\begin{equation}
\label{eq: curvatura constante 2}
\mathbf{R}(u,v)w + \mathbf{R}(u,w)v = K \big( 2\mathbf{g}(v,w)u - \mathbf{g}(u,w)v - \mathbf{g}(u,v)w \big).
\end{equation}

The final step uses the Bianchi identity, which states that
$\mathbf{R}(v,u)w+\mathbf{R}(u,w)v+\mathbf{R}(w,v)u=0$.
Interchanging the roles of $u$ and $v$ in equation (\ref{eq: curvatura constante 2}), we obtain:
\[
\mathbf{R}(v,u)w + \mathbf{R}(v,w)u = K \big( 2\mathbf{g}(u,w)v - \mathbf{g}(v,w)u - \mathbf{g}(u,v)w \big).
\]
By the Bianchi identity and antisymmetry in the first two arguments,
$\mathbf{R}(v,w)u=-\mathbf{R}(w,v)u
=\mathbf{R}(v,u)w+\mathbf{R}(u,w)v$. Substituting gives
\[
\mathbf{R}(v,u)w + \big( \mathbf{R}(v,u)w + \mathbf{R}(u,w)v \big) = 2\mathbf{R}(v,u)w + \mathbf{R}(u,w)v.
\]
Using $\mathbf{R}(v,u)w=-\mathbf{R}(u,v)w$, this is equivalent to
$-2\mathbf{R}(u,v)w+\mathbf{R}(u,w)v$. Therefore,
\begin{equation}
\label{eq: curvatura constante 3}
-2\mathbf{R}(u,v)w + \mathbf{R}(u,w)v = K \big( 2\mathbf{g}(u,w)v - \mathbf{g}(v,w)u - \mathbf{g}(u,v)w \big).
\end{equation}

Finally, subtract equation (\ref{eq: curvatura constante 3}) from equation (\ref{eq: curvatura constante 2}):
\begin{align*}
\big( \mathbf{R}(u,v)w + \mathbf{R}(u,w)v \big) - \big( -2\mathbf{R}(u,v)w + \mathbf{R}(u,w)v \big) &= \\
K \big( 2\mathbf{g}(v,w)u - \mathbf{g}(u,w)v - \mathbf{g}(u,v)w \big) - K \big( 2\mathbf{g}(u,w)v - \mathbf{g}(v,w)u - \mathbf{g}(u,v)w \big).
\end{align*}
Simplifying both sides, the terms involving $\mathbf{R}(u,w)v$ and
$\mathbf{g}(u,v)w$ cancel, leaving
\[
3\mathbf{R}(u,v)w = K \big( 3\mathbf{g}(v,w)u - 3\mathbf{g}(u,w)v \big).
\]
Dividing by $3$ gives the desired formula:
\[
\mathbf{R}(u,v)w = K \big( \mathbf{g}(v,w)u - \mathbf{g}(u,w)v \big).
\]
\end{proof}
\section{Riemannian submanifolds}

Trace conditions and many variational problems are formulated on
a submanifold. To compare intrinsic derivatives with derivatives in $M$,
we must separate the tangential and normal components of the Levi--Civita connection.

Let $(M,\mathbf{g})$ be a Riemannian manifold without boundary and let
$N\subseteq M$ be an embedded submanifold without boundary. The induced metric
is $\mathbf{g}_N=\mathbf{g}\restriction_{TN}$, and
\[
T_xM=T_xN\oplus T_x^\perp N,
\qquad T_x^\perp N=\{v\in T_xM\mid \mathbf{g}(v,w)=0\text{ for every }w\in T_xN\}.
\]
The spaces $T_x^\perp N$ form the normal bundle $T^\perp N$. For a
section $\mathbf{V}$ of $TM|_N$, we write $\mathbf{V}^\top$ and
$\mathbf{V}^\perp$ for its orthogonal projections.

Before differentiating the components of a field, we must check that
the orthogonal decomposition depends smoothly on the base point.

\begin{proposition}[Smoothness of the normal bundle and projections]
\label{prop:haz-normal-proyecciones-suaves}
If $\dim M=n$ and $\dim N=k$, the spaces $T_x^\perp N$ form a
smooth vector subbundle of rank $n-k$ of $TM\restriction_N$.
The orthogonal projections $\Pi^\top$ and $\Pi^\perp$ are smooth
bundle maps. Near each point of $N$, there is a smooth orthonormal
frame $(\mathbf E_1,\ldots,\mathbf E_n)$ of $TM\restriction_N$
such that $\mathbf E_1,\ldots,\mathbf E_k$ are tangent and
$\mathbf E_{k+1},\ldots,\mathbf E_n$ are normal.
\end{proposition}
\begin{proof}
A slice chart provides the tangent frame
$\boldsymbol\partial_1,\ldots,\boldsymbol\partial_k$ and completes
it with the remaining coordinate vector fields of $M$, restricted
to $N$. Apply Gram--Schmidt with the metric $\mathbf g$ in this order.
At each step, subtract the projections onto the vectors already
constructed and divide by the norm of a nonzero vector. The
coefficients are smooth and the denominators positive, so the
resulting frame is smooth. Its first $k$ vectors still span $TN$;
the remaining vectors form an orthonormal basis of the complement.
This basis gives a local trivialization of the normal bundle.
Moreover, for $v\in T_xM$,
\[
\Pi_x^\top v=\displaystyle\sum_{a=1}^k
\mathbf g_x(v,\mathbf E_a(x))\mathbf E_a(x),
\qquad
\Pi_x^\perp v=\displaystyle\sum_{\alpha=k+1}^n
\mathbf g_x(v,\mathbf E_\alpha(x))\mathbf E_\alpha(x).
\]
These expressions are smooth in $(x,v)$ and independent of the
frame, since they characterize the orthogonal projections.
\end{proof}

\begin{lemma}[Ambient connection along a submanifold]
\label{lem:conexion-ambiental-a-lo-largo-subvariedad}
Let $(M,\mathbf{g})$ be a Riemannian manifold without boundary and let
$N\subseteq M$ be an embedded submanifold without boundary. Let
$\iota\colon N\hookrightarrow M$ be the inclusion. There exists a unique connection
$\overline\nabla$ on the restricted bundle $\iota^*TM\to N$ such that, if
$\mathbf{X}\in\mathfrak X(N)$ and $\mathbf{V}\in\Gamma(\iota^*TM)$ admit local extensions
$\widetilde{\mathbf{X}},\widetilde{\mathbf{V}}\in\mathfrak X(M)$, then
\[
\overline\nabla_{\mathbf{X}}\mathbf{V}
=
(\nabla^M_{\widetilde{\mathbf{X}}}\widetilde{\mathbf{V}})\restriction_N.
\]
This expression does not depend on the extensions. Moreover,
\[
\mathbf{X}\mathbf{g}(\mathbf{V},\mathbf{W})
=\mathbf{g}(\overline\nabla_{\mathbf{X}}\mathbf{V},\mathbf{W})
+\mathbf{g}(\mathbf{V},\overline\nabla_{\mathbf{X}}\mathbf{W})
\]
for $\mathbf{V},\mathbf{W}\in\Gamma(\iota^*TM)$ and, if $\mathbf{X},\mathbf{Y}\in\mathfrak X(N)$, then
\[
\overline\nabla_{\mathbf{X}}\mathbf{Y}-\overline\nabla_{\mathbf{Y}}\mathbf{X}=[\mathbf{X},\mathbf{Y}].
\]
\end{lemma}

\begin{proof}
Local extensions exist by
Proposition~\ref{prop:tangentes-extensiones-subvariedad}. In ambient
coordinates $(x^1,\ldots,x^n)$, write
\[
\mathbf V=\displaystyle\sum_{j=1}^n V^j\boldsymbol\partial_j\restriction_N,
\qquad
d\iota(\mathbf X)=\displaystyle\sum_{i=1}^n X^i\boldsymbol\partial_i\restriction_N.
\]
The formula for the ambient connection requires
\[
\overline\nabla_{\mathbf X}\mathbf V
=\displaystyle\sum_{k=1}^n
\left(\mathbf X(V^k)+\displaystyle\sum_{i,j=1}^n
X^iV^j\Gamma^k_{ij}\restriction_N\right)
\boldsymbol\partial_k\restriction_N.
\]
Here $\mathbf X(V^k)$ differentiates a function on $N$. If two
extensions of $V^k$ agree on $N$, their difference vanishes there,
as does its derivative in every tangent direction. The connection
is tensorial in its first argument, so only the value of the
extension of $\mathbf X$ at the point is involved. Thus the
expression is independent of both extensions. Since it is obtained
by restricting an intrinsic ambient connection, expressions in
different charts agree. They define $\overline\nabla$ and also
prove uniqueness. Its rules are
\[
\overline\nabla_{f\mathbf X}\mathbf V
=f\overline\nabla_{\mathbf X}\mathbf V,
\qquad
\overline\nabla_{\mathbf X}(f\mathbf V)
=\mathbf X(f)\mathbf V+f\overline\nabla_{\mathbf X}\mathbf V,
\qquad f\in C^\infty(N),
\]
together with additivity in both arguments.

Metric compatibility follows by differentiating
$\mathbf g(\widetilde{\mathbf V},\widetilde{\mathbf W})$ in the
direction of $\widetilde{\mathbf X}$ and restricting to $N$.
For the torsion identity, choose extensions of tangent fields.
Their bracket restricted to $N$ equals the intrinsic bracket:
in slice coordinates, their normal components vanish on $N$, as do
their tangential derivatives, while the tangential components of the
bracket are precisely those computed on $N$. The torsion-free
identity for $\nabla^M$ therefore gives
$\overline\nabla_{\mathbf X}\mathbf Y-
\overline\nabla_{\mathbf Y}\mathbf X=[\mathbf X,\mathbf Y]$.
\end{proof}

\begin{definition}[Second fundamental form and normal connection]
\label{def:segunda-forma-fundamental-conexion-normal}
Let $(M,\mathbf{g})$ be a Riemannian manifold without boundary and let
$N\subseteq M$ be an embedded submanifold without boundary. For
$\mathbf{X},\mathbf{Y}\in\mathfrak X(N)$ and
$\boldsymbol{\eta}\in\Gamma(T^\perp N)$, define
\[
\mathbf{II}(\mathbf{X},\mathbf{Y})
=(\overline\nabla_{\mathbf{X}}\mathbf{Y})^\perp,
\qquad \nabla_{\mathbf{X}}^\perp\boldsymbol{\eta}
=(\overline\nabla_{\mathbf{X}}\boldsymbol{\eta})^\perp.
\]
For $\boldsymbol{\eta}\in\Gamma(T^\perp N)$, the \textit{Weingarten map} or
\textit{shape operator}
$\mathbf{A}_{\boldsymbol{\eta}}\colon TN\longrightarrow TN$ is the unique
endomorphism
satisfying
\[
\mathbf{g}_N(\mathbf{A}_{\boldsymbol{\eta}}\mathbf{X},\mathbf{Y})
=\mathbf{g}(\mathbf{II}(\mathbf{X},\mathbf{Y}),\boldsymbol{\eta})
\]
for all $\mathbf{X},\mathbf{Y}\in\mathfrak X(N)$.
\end{definition}

\begin{proposition}[Gauss and Weingarten formulas]
\label{prop:gauss-weingarten-subvariedad}
Let $(M,\mathbf{g})$ be a Riemannian manifold without boundary and let
$N\subseteq M$ be an embedded submanifold without boundary. The second fundamental
form $\mathbf{II}$ is a symmetric tensor, and $\nabla^\perp$ is a
metric connection on the normal bundle. We have
\begin{align}
\overline\nabla_{\mathbf{X}}\mathbf{Y}
&=\nabla_{\mathbf{X}}^N\mathbf{Y}+\mathbf{II}(\mathbf{X},\mathbf{Y}),
\label{eq:formula-gauss-subvariedad}\\
\overline\nabla_{\mathbf{X}}\boldsymbol{\eta}
&=-\mathbf{A}_{\boldsymbol{\eta}}\mathbf{X}
+\nabla_{\mathbf{X}}^\perp\boldsymbol{\eta}.
\label{eq:formula-weingarten-subvariedad}
\end{align}
The Weingarten map depends $C^\infty(N)$--linearly on
$\boldsymbol{\eta}$ and is self-adjoint with respect to $\mathbf{g}_N$.
\end{proposition}

\begin{proof}
Lemma~\ref{lem:conexion-ambiental-a-lo-largo-subvariedad} eliminates all
dependence on extensions. Linearity of a connection in its first
argument gives
$\mathbf{II}(f\mathbf{X},\mathbf{Y})=f\mathbf{II}(\mathbf{X},\mathbf{Y})$. In the second argument,
\[
\overline\nabla_{\mathbf{X}}(f\mathbf{Y})
=\mathbf{X}(f)\mathbf{Y}+f\overline\nabla_{\mathbf{X}}\mathbf{Y};
\]
the first summand is tangent, so
$\mathbf{II}(\mathbf{X},f\mathbf{Y})=f\mathbf{II}(\mathbf{X},\mathbf{Y})$.
Thus, $\mathbf{II}$ is tensorial.
Since $\overline\nabla$ is torsion-free on tangent vector fields,
\[
\mathbf{II}(\mathbf{X},\mathbf{Y})-\mathbf{II}(\mathbf{Y},\mathbf{X})
=[\mathbf{X},\mathbf{Y}]^\perp=0.
\]

Define $\mathbf{D}_{\mathbf{X}}\mathbf{Y}
:=(\overline\nabla_{\mathbf{X}}\mathbf{Y})^\top$. Metric compatibility of
$\overline\nabla$ and orthogonality of the decomposition show that $\mathbf{D}$
is metric; the preceding torsion identity shows that it is torsion-free.
By uniqueness of the Levi--Civita connection, $\mathbf{D}=\nabla^N$, and the
tangential--normal decomposition gives the Gauss formula
\eqref{eq:formula-gauss-subvariedad}.

Differentiating the identity
$\mathbf{g}(\boldsymbol{\eta},\mathbf{Y})=0$ gives
\[
\mathbf{g}((\overline\nabla_{\mathbf{X}}\boldsymbol{\eta})^\top,\mathbf{Y})
=-\mathbf{g}(\boldsymbol{\eta},\mathbf{II}(\mathbf{X},\mathbf{Y}))
=-\mathbf{g}_N(\mathbf{A}_{\boldsymbol{\eta}}\mathbf{X},\mathbf{Y}).
\]
Nondegeneracy of $\mathbf{g}_N$ gives
$(\overline\nabla_{\mathbf{X}}\boldsymbol{\eta})^\top
=-\mathbf{A}_{\boldsymbol{\eta}}\mathbf{X}$, proving
\eqref{eq:formula-weingarten-subvariedad}. The connection rules and
metric compatibility of $\nabla^\perp$ follow by projecting those of
$\overline\nabla$. Finally, symmetry of $\mathbf{II}$ implies
\[
\mathbf{g}_N(\mathbf{A}_{\boldsymbol{\eta}}\mathbf{X},\mathbf{Y})
=\mathbf{g}_N(\mathbf{X},\mathbf{A}_{\boldsymbol{\eta}}\mathbf{Y}),
\]
and linearity in $\boldsymbol{\eta}$ follows from the identity defining
$\mathbf{A}_{\boldsymbol{\eta}}$.
\end{proof}

\begin{corollary}[Weingarten map of a Euclidean hypersurface]
\label{cor:weingarten-hipersuperficie-euclidiana}
Let $N\subseteq\mathbb R^{n+1}$ be a smooth oriented hypersurface without
boundary, of dimension $n$, and let
$\boldsymbol{\nu}\colon N\longrightarrow\mathbb S^n$ be its unit normal field.
Then
\[
\mathbf{A}_{\boldsymbol{\nu}}\mathbf{X}=-d\boldsymbol{\nu}(\mathbf{X}),
\qquad \mathbf{X}\in\mathfrak X(N).
\]
In particular, $d\boldsymbol{\nu}$ is self-adjoint with the opposite sign; the eigenvalues
of $\mathbf{A}_{\boldsymbol{\nu}}$ are the principal curvatures, and
$H:=\frac{1}{n}\operatorname{tr}(\mathbf{A}_{\boldsymbol{\nu}})$ is the scalar mean curvature with
this convention.
\end{corollary}

\begin{proof}
Identify $T_{\boldsymbol{\nu}(x)}\mathbb S^n$ with the subspace
$\boldsymbol{\nu}(x)^\perp\subseteq\mathbb R^{n+1}$. Since
$\|\boldsymbol{\nu}\|=1$, we have
$\mathbf{g}(d\boldsymbol{\nu}(\mathbf{X}),\boldsymbol{\nu})
=\frac{1}{2}\mathbf{X}\mathbf{g}(\boldsymbol{\nu},\boldsymbol{\nu})=0$, so
$d\boldsymbol{\nu}(\mathbf{X})$ is
tangent. The Euclidean connection along $N$ satisfies
$\overline\nabla_{\mathbf{X}}\boldsymbol{\nu}
=d\boldsymbol{\nu}(\mathbf{X})$. The Weingarten formula
\eqref{eq:formula-weingarten-subvariedad} and vanishing of the normal
component give
$d\boldsymbol{\nu}(\mathbf{X})=-\mathbf{A}_{\boldsymbol{\nu}}\mathbf{X}$. The remaining assertions are the
corresponding spectral definitions; reality of the eigenvalues
follows from self-adjointness proved in
Proposition~\ref{prop:gauss-weingarten-subvariedad}.
\end{proof}

\begin{definition}[Covariant derivatives of the second fundamental form]
\label{def:derivadas-covariantes-segunda-forma-fundamental}
\index{second fundamental form!covariant derivatives}
Let $(M,\mathbf{g})$ be a Riemannian manifold without boundary and let
$N\subseteq M$ be an embedded submanifold without boundary. The connections
$\nabla^N$ on $TN$ and $\nabla^\perp$ on $T^\perp N$ induce
a connection, denoted by $\nabla^{\mathbf{II}}$, on the bundle
\[
T^*N\otimes T^*N\otimes T^\perp N.
\]
Explicitly, if
$\mathbf{B}\in\Gamma(T^*N\otimes T^*N\otimes T^\perp N)$ and
$\mathbf{X},\mathbf{Y},\mathbf{Z}\in\mathfrak X(N)$, then
\[
(\nabla_{\mathbf{X}}^{\mathbf{II}}\mathbf{B})(\mathbf{Y},\mathbf{Z})
=\nabla_{\mathbf{X}}^\perp(\mathbf{B}(\mathbf{Y},\mathbf{Z}))
-\mathbf{B}(\nabla_{\mathbf{X}}^N\mathbf{Y},\mathbf{Z})
-\mathbf{B}(\mathbf{Y},\nabla_{\mathbf{X}}^N\mathbf{Z}).
\]
Define recursively
\[
(\nabla^{\mathbf{II}})^0\mathbf{II}:=\mathbf{II},
\qquad
(\nabla^{\mathbf{II}})^{m+1}\mathbf{II}
:=\nabla^{\mathbf{II}}
\bigl((\nabla^{\mathbf{II}})^m\mathbf{II}\bigr),
\qquad m\in\mathbb N_0,
\]
where the right-hand side uses the induced connection on each additional
covariant factor. In particular,
\[
(\nabla^{\mathbf{II}})^m\mathbf{II}
\in
\Gamma\bigl((T^*N)^{\otimes(m+2)}\otimes T^\perp N\bigr).
\]
The dual metric of $\mathbf{g}_N$ on the factors $T^*N$ and the restriction of $\mathbf{g}$ to the
normal bundle induce the pointwise norm
$|(\nabla^{\mathbf{II}})^m\mathbf{II}|_{\mathbf{g}}$. We say that this derivative is
uniformly bounded if
\[
\displaystyle\sup_{x\in N}
\bigl|(\nabla^{\mathbf{II}})^m\mathbf{II}(x)\bigr|_{\mathbf{g}}<\infty.
\]
\end{definition}

\begin{definition}[Normal curvature]
Let $(M,\mathbf{g})$ be a Riemannian manifold without boundary and let
$N\subseteq M$ be an embedded submanifold without boundary. The curvature of the
normal connection is the field $\mathbf{R}^\perp$ defined by
\[
\mathbf{R}^\perp(\mathbf{X},\mathbf{Y})\boldsymbol{\eta}
:=
\nabla_{\mathbf{X}}^\perp\nabla_{\mathbf{Y}}^\perp\boldsymbol{\eta}
-\nabla_{\mathbf{Y}}^\perp\nabla_{\mathbf{X}}^\perp\boldsymbol{\eta}
-\nabla_{[\mathbf{X},\mathbf{Y}]}^\perp\boldsymbol{\eta}.
\]
\end{definition}

\begin{theorem}[Gauss, Codazzi, and Ricci equations]
\label{teo:gauss-codazzi-ricci-subvariedad}
Let $(M,\mathbf{g})$ be a Riemannian manifold without boundary and let
$N\subseteq M$ be an embedded submanifold without boundary. For tangent vector fields
$\mathbf{X},\mathbf{Y},\mathbf{Z},\mathbf{W}$ and normal fields
$\boldsymbol{\eta},\boldsymbol{\zeta}$, we have
\begin{align}
\mathbf{g}(\mathbf{R}^N(\mathbf{X},\mathbf{Y})\mathbf{Z},\mathbf{W})
={}&\mathbf{g}(\mathbf{R}^M(\mathbf{X},\mathbf{Y})\mathbf{Z},\mathbf{W})
+\mathbf{g}(\mathbf{II}(\mathbf{X},\mathbf{W}),
\mathbf{II}(\mathbf{Y},\mathbf{Z}))
-\mathbf{g}(\mathbf{II}(\mathbf{X},\mathbf{Z}),
\mathbf{II}(\mathbf{Y},\mathbf{W})),\\
(\mathbf{R}^M(\mathbf{X},\mathbf{Y})\mathbf{Z})^\perp
={}&(\nabla_{\mathbf{X}}^{\mathbf{II}}\mathbf{II})(\mathbf{Y},\mathbf{Z})
-(\nabla_{\mathbf{Y}}^{\mathbf{II}}\mathbf{II})(\mathbf{X},\mathbf{Z}),\\
\mathbf{g}(\mathbf{R}^\perp(\mathbf{X},\mathbf{Y})\boldsymbol{\eta},
\boldsymbol{\zeta})
={}&\mathbf{g}(\mathbf{R}^M(\mathbf{X},\mathbf{Y})\boldsymbol{\eta},
\boldsymbol{\zeta})
+\mathbf{g}_N([\mathbf{A}_{\boldsymbol{\eta}},\mathbf{A}_{\boldsymbol{\zeta}}]
\mathbf{X},\mathbf{Y}).
\end{align}
We use the convention
$\mathbf{R}(\mathbf{X},\mathbf{Y})\mathbf{Z}
=\nabla_{\mathbf{X}}\nabla_{\mathbf{Y}}\mathbf{Z}
-\nabla_{\mathbf{Y}}\nabla_{\mathbf{X}}\mathbf{Z}
-\nabla_{[\mathbf{X},\mathbf{Y}]}\mathbf{Z}$ from
Definition~\ref{def:variedades-riemannianas-variedad-5}.
\end{theorem}

\begin{proof}
We may compute ambient curvature using fields along $N$ by
Lemma~\ref{lem:conexion-ambiental-a-lo-largo-subvariedad}.
Apply Gauss and Weingarten twice. For tangent fields,
\[
\begin{aligned}
\overline\nabla_{\mathbf X}\overline\nabla_{\mathbf Y}\mathbf Z
={}&\nabla^N_{\mathbf X}\nabla^N_{\mathbf Y}\mathbf Z
-\mathbf A_{\mathbf{II}(\mathbf Y,\mathbf Z)}\mathbf X\\
&+\mathbf{II}(\mathbf X,\nabla^N_{\mathbf Y}\mathbf Z)
+\nabla^\perp_{\mathbf X}(\mathbf{II}(\mathbf Y,\mathbf Z)).
\end{aligned}
\]
The first two terms are tangent and the last two normal. Subtract
the expression with $\mathbf X,\mathbf Y$ interchanged, and subtract
$\overline\nabla_{[\mathbf X,\mathbf Y]}\mathbf Z$.
The tangential component is
\[
(\mathbf R^M(\mathbf X,\mathbf Y)\mathbf Z)^\top
=\mathbf R^N(\mathbf X,\mathbf Y)\mathbf Z
-\mathbf A_{\mathbf{II}(\mathbf Y,\mathbf Z)}\mathbf X
+\mathbf A_{\mathbf{II}(\mathbf X,\mathbf Z)}\mathbf Y.
\]
Taking its inner product with $\mathbf W$ and using the definition
of the shape operator yields the Gauss equation with the stated signs.

The normal component of the same calculation is
\[
\begin{aligned}
(\mathbf R^M(\mathbf X,\mathbf Y)\mathbf Z)^\perp
={}&\nabla^\perp_{\mathbf X}(\mathbf{II}(\mathbf Y,\mathbf Z))
-\nabla^\perp_{\mathbf Y}(\mathbf{II}(\mathbf X,\mathbf Z))\\
&+\mathbf{II}(\mathbf X,\nabla^N_{\mathbf Y}\mathbf Z)
-\mathbf{II}(\mathbf Y,\nabla^N_{\mathbf X}\mathbf Z)
-\mathbf{II}([\mathbf X,\mathbf Y],\mathbf Z).
\end{aligned}
\]
Expanding the definition of $\nabla^{\mathbf{II}}$, the difference
$(\nabla^{\mathbf{II}}_{\mathbf X}\mathbf{II})(\mathbf Y,\mathbf Z)
-(\nabla^{\mathbf{II}}_{\mathbf Y}\mathbf{II})(\mathbf X,\mathbf Z)$
has the same first four terms and the additional term
\[
-\mathbf{II}(\nabla^N_{\mathbf X}\mathbf Y
-\nabla^N_{\mathbf Y}\mathbf X,\mathbf Z)
=-\mathbf{II}([\mathbf X,\mathbf Y],\mathbf Z).
\]
This equality uses the absence of torsion of $\nabla^N$ and proves
Codazzi.

For a normal field $\boldsymbol\eta$, the normal component of the
second derivative is
\[
(\overline\nabla_{\mathbf X}\overline\nabla_{\mathbf Y}
\boldsymbol\eta)^\perp
=-\mathbf{II}(\mathbf X,\mathbf A_{\boldsymbol\eta}\mathbf Y)
+\nabla^\perp_{\mathbf X}\nabla^\perp_{\mathbf Y}\boldsymbol\eta.
\]
Antisymmetrizing and subtracting the bracket term gives
\[
(\mathbf R^M(\mathbf X,\mathbf Y)\boldsymbol\eta)^\perp
=\mathbf R^\perp(\mathbf X,\mathbf Y)\boldsymbol\eta
-\mathbf{II}(\mathbf X,\mathbf A_{\boldsymbol\eta}\mathbf Y)
+\mathbf{II}(\mathbf Y,\mathbf A_{\boldsymbol\eta}\mathbf X).
\]
Finally, self-adjointness of both shape operators gives
\[
\begin{aligned}
&\mathbf g(\mathbf{II}(\mathbf X,\mathbf A_{\boldsymbol\eta}\mathbf Y),
\boldsymbol\zeta)
-\mathbf g(\mathbf{II}(\mathbf Y,\mathbf A_{\boldsymbol\eta}\mathbf X),
\boldsymbol\zeta)\\
&\quad=\mathbf g_N(\mathbf A_{\boldsymbol\zeta}\mathbf X,
\mathbf A_{\boldsymbol\eta}\mathbf Y)
-\mathbf g_N(\mathbf A_{\boldsymbol\zeta}\mathbf Y,
\mathbf A_{\boldsymbol\eta}\mathbf X)\\
&\quad=\mathbf g_N\bigl((\mathbf A_{\boldsymbol\eta}\mathbf A_{\boldsymbol\zeta}
-\mathbf A_{\boldsymbol\zeta}\mathbf A_{\boldsymbol\eta})\mathbf X,
\mathbf Y\bigr).
\end{aligned}
\]
Solving for the inner product involving $\mathbf R^\perp$ proves Ricci.
\end{proof}

The Gauss equation shows how intrinsic curvature records both
ambient curvature and the bending of the submanifold. For
orthonormal $X,Y\in T_xN$, setting $\mathbf Z=\mathbf Y$ and
$\mathbf W=\mathbf X$ gives
\[
K_N(\operatorname{span}\{X,Y\})
=K_M(\operatorname{span}\{X,Y\})
+\mathbf g_x(\mathbf{II}(X,X),\mathbf{II}(Y,Y))
-|\mathbf{II}(X,Y)|_{\mathbf g_x}^2.
\]
Here $K_M$ is evaluated on the same plane, viewed as a subspace of
$T_xM$.

\begin{definition}[Mean curvature and totally geodesic submanifolds]
\label{def:curvatura-media-totalmente-geodesica}
For $k=\dim N\geq1$, the \textbf{mean curvature vector} is
\[
\mathbf H:=\frac1k\displaystyle\sum_{a=1}^k
\mathbf{II}(\mathbf E_a,\mathbf E_a),
\]
where $(\mathbf E_a)_{a=1}^k$ is any tangent orthonormal frame.
The submanifold is \textbf{totally geodesic} if $\mathbf{II}=0$.
\end{definition}
The definition of $\mathbf H$ is independent of the frame: if
$\mathbf F_a=\displaystyle\sum_{b=1}^k O_a^{,b}\mathbf E_b$ with
$O$ orthogonal, then
$\displaystyle\sum_{a=1}^k O_a^{,b}O_a^{,c}=\delta_{bc}$, and
bilinearity of $\mathbf{II}$ gives the same sum. For a hypersurface
with unit normal $\boldsymbol\nu$, we have
$\mathbf H=H\boldsymbol\nu$, with the preceding scalar convention.

\begin{proposition}[Geodesic characterization]
\label{prop:totalmente-geodesica-caracterizacion}
The submanifold $N$ is totally geodesic if and only if every geodesic
of $(N,\mathbf g_N)$, with its affine parameter, is also a geodesic
of $(M,\mathbf g)$.
\end{proposition}
\begin{proof}
The Gauss formula along a curve $\gamma$ in $N$ gives
\[
D_t^M\dot\gamma=D_t^N\dot\gamma
+\mathbf{II}(\dot\gamma,\dot\gamma).
\]
This follows by locally extending the velocity or by applying the
coordinate formula; both sides depend only on the curve.
If $\mathbf{II}=0$, a geodesic of $N$ has zero ambient acceleration.
Conversely, for every $v\in T_xN$ there is a local geodesic of $N$
with initial velocity $v$. If it is also an ambient geodesic, the
identity gives $\mathbf{II}_x(v,v)=0$. Symmetry and polarization yield
\[
2\mathbf{II}_x(v,w)=\mathbf{II}_x(v+w,v+w)
-\mathbf{II}_x(v,v)-\mathbf{II}_x(w,w)=0.
\]
Thus $\mathbf{II}=0$.
\end{proof}

For example, on the sphere of radius $R$ in $\mathbb R^{k+1}$,
with outward normal $\boldsymbol\nu(x)=x/R$, we have
$\mathbf A_{\boldsymbol\nu}=-R^{-1}\operatorname{Id}$ and
$\mathbf{II}(X,Y)=-R^{-1}\mathbf g(X,Y)\boldsymbol\nu$.
Thus $\mathbf H=-R^{-1}\boldsymbol\nu$ and, for $k\geq2$,
the Gauss formula gives $K_N=R^{-2}$. This explicitly fixes the
signs in our convention.

If the submanifold has boundary, the preceding tensor identities
remain valid up to it: prove them in the interior and extend them
by continuity in smooth boundary charts. To treat $\partial M$
when the ambient manifold itself has boundary, extend the metric
locally to the double using
Corollary~\ref{cor:extension-suave-al-doble}. The values and
derivatives of two extensions agree on $M$ by continuity from its
interior; hence the connections and fundamental forms obtained on
$\partial M$ are independent of that choice.

\section{Tubular neighborhoods and Fermi coordinates}

The exponential map applied to the normal bundle turns the preceding tangential--normal decomposition into coordinates adapted to $N$.

We may adapt the exponential map to obtain information about the normal bundle as follows, defining what is known as the \textit{normal exponential map}.

\begin{definition}\label{def:variedades-riemannianas-variedad-3}\index{normal exponential map}
Let $(M,\mathbf{g})$ be a Riemannian manifold without boundary and let
$N\subseteq M$ be an embedded submanifold without boundary. Define the normal
exponential map
$\operatorname{exp}^{\perp}\colon \mathcal{E}_{N}\longrightarrow M$, with
$\mathcal{E}_{N}:=\mathcal{E}\cap T^{\perp}N$, by
$\exp^{\perp}=\operatorname{exp}|_{\mathcal{E}_{N}}$.
\end{definition}

We will be interested in ``normal neighborhoods around $N$,'' which we call \textit{tubular neighborhoods}, as in the following figure:

\begin{figura}[H]
    \includegraphics[width=\linewidth]{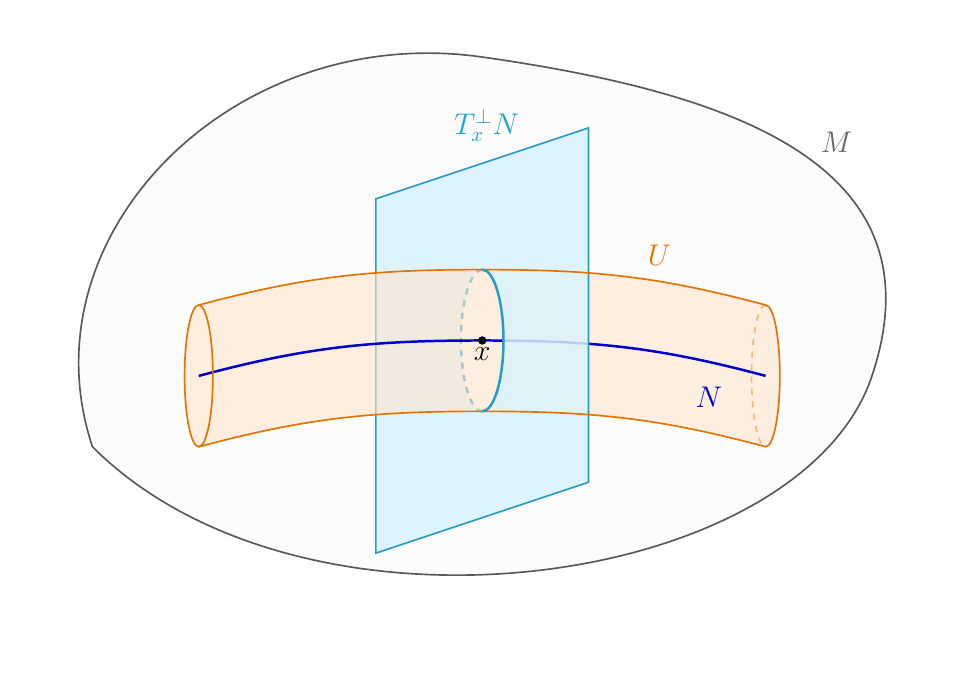}
    \caption{Tubular neighborhood around the embedded submanifold $N$.}
    \label{fig:vecindad-tubular}
\end{figura}

\begin{definition}\label{def:variedades-riemannianas-variedad-4}\index{tubular neighborhood}
Let $(M,\mathbf{g})$ be a Riemannian manifold without boundary and let
$N\subseteq M$ be an embedded submanifold without boundary. We say that
$U\subseteq M$ is a \textit{normal neighborhood} of $N$ in $M$ if there exists an
open set $V\subseteq \mathcal{E}_{N}$ such that:
\begin{enumerate}
 \item $V$ intersects each fiber $T^{\perp}_{x}N$ in a set star-shaped with respect to $0$,
 \item the restriction
 \[
 \exp^{\perp}\restriction_{V}\colon V\longrightarrow U
 \]
 is a diffeomorphism.
\end{enumerate}

We say that $U$ is a \textit{tubular neighborhood} if it is a normal neighborhood and, in addition, $V$ can be chosen in the form
\[
V:=\{(x,v)\in T^{\perp}N \mid |v|_{\mathbf{g}}<\delta(x)\},
\]
for some continuous function $\delta\colon N\longrightarrow (0,\infty)$.

If $V$ can be chosen with $\delta(x)\equiv \varepsilon$, we say that $U$ is a \textit{uniform tubular neighborhood} of radius $\varepsilon$, or an $\varepsilon$-tubular neighborhood.
\end{definition}

Every embedded submanifold admits a tubular neighborhood, and if the submanifold is compact, this tubular neighborhood can be chosen to be uniform:

\begin{theorem}[Existence of tubular neighborhoods]\label{teo:variedades-riemannianas-existencia-de-vecindades-tubulares}\index{existence of tubular neighborhoods}
Let $(M,\mathbf{g})$ be a Riemannian manifold without boundary and let
$N\subseteq M$ be an embedded submanifold without boundary. Then $N$ has
a tubular neighborhood in $M$. If $N$ is compact, it has a uniform tubular
neighborhood.
\end{theorem}

\begin{proof}
Let \(N_0:=\{(x,0)\mid x\in N\}\subseteq T^\perp N\). At a point of the
zero section, there is a canonical decomposition
\begin{equation}
T_{(x,0)}(T^\perp N)
=T_{(x,0)}N_0\oplus T_{(x,0)}(T_x^\perp N)
\cong T_xN\oplus T_x^\perp N.
\label{eq:tangente-haz-normal-cero}
\end{equation}
The restriction of \(\exp^\perp\) to \(N_0\) is the inclusion \(N\hookrightarrow M\);
its differential maps the first summand onto \(T_xN\). On the fiber
\(T_x^\perp N\), the normal exponential map equals \(\exp_x\), and
\(d(\exp_x)_0=I\) by Proposition~\ref{propiedades mapeo exponencial}. Therefore, \(d(\exp^\perp)_{(x,0)}\) is the isomorphism
\[
T_xN\oplus T_x^\perp N\longrightarrow T_xM,
\qquad (X,V)\longmapsto X+V.
\]
The inverse function theorem, Theorem~\ref{teo: funcion inversa variedades},
shows that \(\exp^\perp\) is a local diffeomorphism near each point of
\(N_0\).

To obtain a single neighborhood of the entire zero section, we must prevent
normal geodesics with different initial points from intersecting. Denote
the ambient distance restricted to \(N\times N\) by \(d_{\mathbf{g}}\) and set
\[
V_\delta(x):=
\{(y,v)\in T^\perp N\mid d_{\mathbf{g}}(x,y)<\delta,\ |v|_{\mathbf{g}}<\delta\},
\qquad 0<\delta\leq1.
\]
For each \(x\), the preceding local result allows us to choose \(\delta>0\) so
that \(\exp^\perp\) is a diffeomorphism from \(V_\delta(x)\) onto its image.
Define
\begin{equation}
\rho(x):=\displaystyle\sup\{0<\delta\leq1\mid
\exp^\perp\restriction_{V_\delta(x)}
\text{ is a diffeomorphism onto its image}\}.
\label{eq:radio-local-tubular}
\end{equation}
Then \(\rho(x)>0\). The property defining the supremum is monotone with respect
to \(\delta\). Hence, given two points of \(V_{\rho(x)}(x)\), there exist
\(\delta<\rho(x)\) containing both and an \(\delta'\), with
\(\delta<\delta'<\rho(x)\), for which the restriction is a diffeomorphism.
Thus, \(\exp^\perp\) is injective and a local diffeomorphism on
\(V_{\rho(x)}(x)\), and therefore is a diffeomorphism from this set onto its
image.

The function \(\rho\) is continuous. If \(d_{\mathbf{g}}(x,x')<\rho(x)\), the
triangle inequality gives
\[
V_{\rho(x)-d_{\mathbf{g}}(x,x')}(x')\subseteq V_{\rho(x)}(x),
\]
and hence
\(\rho(x')\geq\rho(x)-d_{\mathbf{g}}(x,x')\). If \(d_{\mathbf{g}}(x,x')\geq\rho(x)\), this same
inequality is trivial. Interchanging \(x\) and \(x'\) gives
\begin{equation}
|\rho(x)-\rho(x')|\leq d_{\mathbf{g}}(x,x').
\label{eq:rho-lipschitz-tubular}
\end{equation}
The ambient distance induces the topology of the embedded submanifold,
so \(\rho\) is continuous.

Consider the open set
\begin{equation}
V:=\{(x,v)\in T^\perp N\mid |v|_{\mathbf{g}}<\tfrac12\rho(x)\}.
\label{eq:dominio-tubular-global}
\end{equation}
Let us prove that \(\exp^\perp\) is injective on \(V\). Suppose that
\(\exp_x(v)=\exp_{x'}(v')\) and, without loss of generality,
\(\rho(x')\leq\rho(x)\). The two radial geodesics, the second traversed in
reverse, form an admissible curve from \(x\) to \(x'\) of length
\(|v|_{\mathbf{g}}+|v'|_{\mathbf{g}}\). Consequently,
\[
d_{\mathbf{g}}(x,x')\leq |v|_{\mathbf{g}}+|v'|_{\mathbf{g}}
<\tfrac12\rho(x)+\tfrac12\rho(x')
\leq\rho(x).
\]
Both \((x,v)\) and \((x',v')\) then belong to
\(V_{\rho(x)}(x)\), where \(\exp^\perp\) is injective; the two points are
equal.

Every point of \(V\) belongs to an open set on which \(\exp^\perp\) is a
local diffeomorphism. Consequently,
\(\exp^\perp\restriction_V\) is an injective local diffeomorphism, hence a
diffeomorphism onto the open set \(U:=\exp^\perp(V)\). The formulas
\eqref{eq:dominio-tubular-global} and continuity of \(\frac{\rho}{2}\) show that
\(U\) is a tubular neighborhood. If \(N\) is compact, \(\frac{\rho}{2}\) attains a
positive minimum \(\varepsilon\), and restriction to
\(\{|v|_{\mathbf{g}}<\varepsilon\}\) gives the uniform tubular neighborhood.
\end{proof}
\subsection{Fermi coordinates}

The tubular neighborhood identifies a neighborhood of \(N\) with a neighborhood of the zero
section of its normal bundle. Local trivializations of this bundle give coordinates
adapted simultaneously to the tangential and normal directions.

\begin{definition}[Fermi coordinates]
\label{def:coordenadas-fermi-locales}
\index{Fermi coordinates}
Let \((M,\mathbf{g})\) be a Riemannian manifold without boundary of dimension
\(n\), let \(N\subseteq M\) be an embedded submanifold without boundary of dimension
\(k\), and let \(U=\exp^\perp(V)\) be a normal neighborhood of \(N\). Let
\(\theta=(y^1,\dots,y^k)\colon W\longrightarrow\theta(W)\subseteq\mathbb R^k\) be a chart of \(N\) and let
\((\mathbf{E}_1,\dots,\mathbf{E}_c)\), \(c=n-k\), be a local orthonormal frame of \(T^\perp N\)
over \(W\). On any open product \(W'\times Z\subseteq
\theta(W)\times\mathbb R^c\) on which the following expression takes values in
\(V\), define
\begin{equation}
\Phi(u,z):=
\exp_{\,\theta^{-1}(u)}^\perp
\left(\displaystyle\sum_{\alpha=1}^c z^\alpha
\mathbf{E}_\alpha\big|_{\theta^{-1}(u)}\right).
\label{eq:parametrizacion-fermi}
\end{equation}
The inverse chart of \(\Phi\) is called a \textit{Fermi chart}; its
coordinates are denoted by
\((y^1,\dots,y^k,z^1,\dots,z^c)\).
\end{definition}

The parametrization \eqref{eq:parametrizacion-fermi} is a diffeomorphism onto
its image: it is the composition of the local trivialization of the normal bundle with
\(\exp^\perp\restriction_V\). Thus the definition presupposes no
additional property of the coordinates.

\begin{proposition}[Local properties of Fermi coordinates]
\label{prop:propiedades-coordenadas-fermi-locales}
Let \((M,\mathbf{g})\) be a Riemannian manifold without boundary of dimension
\(n\), let \(N\subseteq M\) be an embedded submanifold without boundary of dimension
\(k\), and let \((y^1,\ldots,y^k,z^1,\ldots,z^c)\) be a Fermi chart associated with the
local normal frame \((\mathbf{E}_1,\ldots,\mathbf{E}_c)\). With indices
\(1\leq a,b\leq k\) and
\(1\leq\alpha,\beta,\gamma\leq c\), the following hold:
\begin{enumerate}[label=(\alph*)]
\item \(N\) is described by \(z^1=\dots=z^c=0\).
\item On \(N\),
\begin{equation}
g_{a,k+\alpha}=0,
\qquad
g_{k+\alpha,k+\beta}=\delta_{\alpha\beta}.
\label{eq:bloques-metrica-fermi}
\end{equation}
\item For \(q\in W\) and \(v=\displaystyle\sum_{\alpha=1}^c z^\alpha \mathbf{E}_\alpha|_q\), the geodesic with initial
data \((q,v)\) has the expression
\begin{equation}
t\longmapsto
(y^1(q),\dots,y^k(q),tz^1,\dots,tz^c).
\label{eq:geodesica-normal-fermi}
\end{equation}
for those values of $t$ for which $(q,tv)$ belongs to the domain of
the chosen tubular parametrization.
\item On \(N\),
\begin{equation}
\Gamma^i_{k+\alpha,k+\beta}=0
\quad
(1\leq i\leq n,\ 1\leq\alpha,\beta\leq c).
\label{eq:christoffel-normal-fermi}
\end{equation}
\item On \(N\),
\begin{equation}
\partial_{z^\alpha}g_{k+\beta,k+\gamma}=0.
\label{eq:derivada-metrica-normal-fermi}
\end{equation}
\end{enumerate}
\end{proposition}

The local identities agree with Proposition~5.26 in \cite{LeeR}.

\begin{proof}
Part (a) follows from \(\exp_q(0)=q\) and injectivity of the normal exponential
map on the chosen neighborhood. At \((q,0)\), the differential of \(\Phi\) sends
\(\boldsymbol{\partial}_{y^a}\) to a vector in \(T_qN\) and
\(\boldsymbol{\partial}_{z^\alpha}\) to
\(\mathbf{E}_\alpha|_q\in T_q^\perp N\). Orthogonality and
orthonormality of the frame prove \eqref{eq:bloques-metrica-fermi}.

The rescaling lemma for geodesics,
Lemma~\ref{lema: reescalamiento geodesicas}, gives
\(\exp_q(tv)=\gamma_v(t)\) and proves
\eqref{eq:geodesica-normal-fermi}. Substituting this curve into the geodesic
equation and evaluating at \(t=0\) gives, for each \(i\),
\[
\displaystyle\sum_{\alpha,\beta=1}^c\Gamma^i_{k+\alpha,k+\beta}(q,0)z^\alpha z^\beta=0
\quad\text{for every }z\in\mathbb R^c.
\]
The Levi--Civita connection is torsion-free, so the two lower indices
of \(\Gamma\) are symmetric. Polarizing the preceding quadratic form
proves \eqref{eq:christoffel-normal-fermi}. Finally, metric
compatibility in coordinates, Corollary~\ref{compatibilidad levi civita metrica en coordenadas},
gives
\[
\partial_{z^\alpha}g_{k+\beta,k+\gamma}
=
\displaystyle\sum_{m=1}^n\bigl(\Gamma^m_{k+\alpha,k+\beta}g_{m,k+\gamma}
+\Gamma^m_{k+\alpha,k+\gamma}g_{k+\beta,m}\bigr),
\]
and both terms vanish by
\eqref{eq:christoffel-normal-fermi}. This proves
\eqref{eq:derivada-metrica-normal-fermi}.
\end{proof}

\begin{proposition}[Boundary normal coordinates]
\label{obs:coordenadas-normales-frontera}
Let \((M,\mathbf{g})\) be a smooth Riemannian manifold of dimension \(n\)
with nonempty boundary and let
\(p\in\partial M\). There exist a chart
\(\theta=(y^1,\dots,y^{n-1})\colon W\longrightarrow\theta(W)\) of \(\partial M\) around \(p\), an
\(\varepsilon>0\), and a parametrization whose inverse is a chart of \(M\)
\begin{equation}
\Phi(u,t)=\exp_{\theta^{-1}(u)}
\bigl(t\boldsymbol{\nu}_{\theta^{-1}(u)}\bigr),
\qquad 0\leq t<\varepsilon,
\label{eq:coordenadas-normales-frontera}
\end{equation}
where \(\boldsymbol{\nu}\) is the inward unit normal. Throughout the chart,
\begin{equation}
\mathbf{g}=dt^2+\mathbf{g}_t,
\qquad g_{tt}=1,
\qquad g_{at}=0.
\label{eq:forma-semigeodesica-metrica}
\end{equation}
If
\(\mathbf{II}^{\boldsymbol{\nu}}(\mathbf{X},\mathbf{Y})
:=\mathbf{g}(\mathbf{II}(\mathbf{X},\mathbf{Y}),\boldsymbol{\nu})\) and we use
the convention
\(\nabla_{\mathbf{X}}\boldsymbol{\nu}
=-\mathbf{A}_{\boldsymbol{\nu}}\mathbf{X}\), then
\begin{equation}
\left.\partial_tg_{ab}\right|_{t=0}
=-2\mathbf{II}_{ab}^{\boldsymbol{\nu}}.
\label{eq:derivada-normal-metrica-II}
\end{equation}
\end{proposition}

\begin{proof}
By Corollary~\ref{cor:extension-suave-al-doble}, extend $\mathbf g$
to a neighborhood of $M$ in its double. The inward unit normal of
$\partial M$ is smooth: in a chart of the double adapted to the copy $M$, described by
$r\geq0$, it is the metric dual of $dr$, normalized to length one,
and satisfies $dr(\boldsymbol\nu)>0$. The local expressions agree
on overlaps because they select the same inward normal direction.

Apply the tubular neighborhood theorem to $\partial M$ in this
extension. It gives an injective normal parametrization near the
zero section. Choose a chart $\theta$ of $\partial M$ with domain
$W$ whose closure is compact in a smaller boundary neighborhood.
Compactness gives a common $\varepsilon>0$ such that
$(x,t\boldsymbol\nu_x)$ belongs to the tubular domain for every
$x\in W$ and $|t|<\varepsilon$. Shrink $W$ around $p$ so that the
image lies in a chart of the double adapted to $M$, where $M$ is given by $r\geq0$. There,
\[
\left.\frac{d}{dt}\right|_{t=0}
 r(\exp_x(t\boldsymbol\nu_x))=dr_x(\boldsymbol\nu_x)>0.
\]
Continuity and compactness allow us to reduce $\varepsilon$ so that
this derivative remains positive for $x\in\overline W$ and
$|t|<\varepsilon$. Since the initial value of $r$ is zero,
$r(\exp_x(t\boldsymbol\nu_x))$ has the sign of $t$. Restricting the
tubular parametrization to $t\geq0$ therefore gives a diffeomorphism
onto a relative neighborhood of $p$ in $M$. This proves the existence
of \eqref{eq:coordenadas-normales-frontera}.

Let \(\mathbf{T}=\boldsymbol{\partial}_t\) and let
\(\mathbf{J}_a=\boldsymbol{\partial}_{y^a}\) be the variational field obtained
by varying the initial point of the normal geodesics. Each curve
\(t\mapsto\Phi(y,t)\) is a unit-speed geodesic, so
\(\mathbf{g}(\mathbf{T},\mathbf{T})=1\). The
coordinate vector fields commute and the connection is torsion-free; therefore,
\[
\frac{d}{dt}\mathbf{g}(\mathbf{J}_a,\mathbf{T})
=\mathbf{g}(\nabla_{\mathbf{T}}\mathbf{J}_a,\mathbf{T})
+\mathbf{g}(\mathbf{J}_a,\nabla_{\mathbf{T}}\mathbf{T})
=\mathbf{g}(\nabla_{\mathbf{J}_a}\mathbf{T},\mathbf{T})
=\frac{1}{2}\mathbf{J}_a\mathbf{g}(\mathbf{T},\mathbf{T})=0.
\]
At \(t=0\), \(\mathbf{J}_a\) is tangent to \(\partial M\) and
\(\mathbf{T}=\boldsymbol{\nu}\), so
\(\mathbf{g}(\mathbf{J}_a,\mathbf{T})=0\). This proves
\eqref{eq:forma-semigeodesica-metrica} throughout the collar.

Using \([\mathbf{T},\mathbf{J}_a]=0\) and metric compatibility again,
\[
\partial_tg_{ab}
=\mathbf{g}(\nabla_{\mathbf{T}}\mathbf{J}_a,\mathbf{J}_b)
+\mathbf{g}(\mathbf{J}_a,\nabla_{\mathbf{T}}\mathbf{J}_b)
=\mathbf{g}(\nabla_{\mathbf{J}_a}\mathbf{T},\mathbf{J}_b)
+\mathbf{g}(\mathbf{J}_a,\nabla_{\mathbf{J}_b}\mathbf{T}).
\]
At \(t=0\), the Weingarten identity
\(\nabla_{\mathbf{J}_a}\boldsymbol{\nu}
=-\mathbf{A}_{\boldsymbol{\nu}}\mathbf{J}_a\) gives
\eqref{eq:derivada-normal-metrica-II}.
Indeed, at $t=0$ the two terms are $-\mathbf{II}^{\boldsymbol\nu}_{ab}$ and $-\mathbf{II}^{\boldsymbol\nu}_{ba}$; they agree by symmetry of the second fundamental form.
\end{proof}
\section{Jacobi fields and comparison theory}

Jacobi fields describe the differential of the exponential map and therefore
the infinitesimal separation of nearby geodesics.

\begin{proposition}[Commutation in a variation]
\label{prop:conmutacion-variacion-jacobi}
Let \((M,\mathbf g)\) be a Riemannian manifold without boundary, let
\(\alpha\) be a smooth variation, and let \(\mathbf Z\) be a smooth
section of \(\alpha^*TM\). If
\(\mathbf S=d\alpha(\boldsymbol{\partial}_s)\) and
\(\mathbf T=d\alpha(\boldsymbol{\partial}_t)\), then
\[
 D_s\mathbf T=D_t\mathbf S,\qquad
 (D_sD_t-D_tD_s)\mathbf Z=\mathbf R(\mathbf S,\mathbf T)\mathbf Z.
\]
\end{proposition}
\begin{proof}
The first identity is
Lemma~\ref{lem:conmutacion-velocidades-variacion}. For the second,
work in a chart on \(M\), with \(n=\dim M\). Write
\(\mathbf Z=\displaystyle\sum_{j=1}^n Z^j\boldsymbol{\partial}_j|_\alpha\)
and define the matrices
\[
 (A_s)^k_j:=\displaystyle\sum_{i=1}^n
 \Gamma^k_{ij}(\alpha)S^i,\qquad
 (A_t)^k_j:=\displaystyle\sum_{i=1}^n
 \Gamma^k_{ij}(\alpha)T^i,
 \qquad 1\leq j,k\leq n.
\]
In components, \(D_sZ=\partial_sZ+A_sZ\) and
\(D_tZ=\partial_tZ+A_tZ\). Expanding both products, all terms
containing first derivatives of \(Z\) cancel, leaving
\[
 (D_sD_t-D_tD_s)Z
 =(\partial_sA_t-\partial_tA_s+A_sA_t-A_tA_s)Z.
\]
Since \(\partial_sT^i=\partial_tS^i\), component \(k\) of this vector is
\[
 \displaystyle\sum_{i,\ell,j=1}^n
 C^k_{i\ell j}(\alpha)S^iT^\ell Z^j,
\]
where
\[
 C^k_{i\ell j}
 :=\partial_i\Gamma^k_{\ell j}-\partial_\ell\Gamma^k_{ij}
 +\displaystyle\sum_{m=1}^n
 \bigl(\Gamma^k_{im}\Gamma^m_{\ell j}
       -\Gamma^k_{\ell m}\Gamma^m_{ij}\bigr).
\]
The definition of curvature applied to coordinate vector fields gives
\[
 \mathbf R(\boldsymbol{\partial}_i,\boldsymbol{\partial}_\ell)
 \boldsymbol{\partial}_j
 =\displaystyle\sum_{k=1}^n C^k_{i\ell j}\boldsymbol{\partial}_k.
\]
Thus the preceding expression is
\(\mathbf R(\mathbf S,\mathbf T)\mathbf Z\).
The calculation does not require \(\alpha\) to be an immersion.
\end{proof}

\begin{theorem}[Jacobi equation]
\label{teo:ecuacion-jacobi-variacion-geodesica}
\index{Jacobi equation}\index{Jacobi field}
Let $(M,\mathbf{g})$ be a Riemannian manifold without boundary. The variational
field $\mathbf{J}$ of a variation through geodesics satisfies
\begin{equation}
 \frac{D^2\mathbf{J}}{dt^2}+\mathbf{R}(\mathbf{J},\dot\gamma)\dot\gamma=0.
 \label{eq:ecuacion-jacobi-convencion-curvatura}
\end{equation}
A smooth field along a geodesic satisfying
\eqref{eq:ecuacion-jacobi-convencion-curvatura} is called a
\textbf{Jacobi field}.
\end{theorem}
\begin{proof}
For each $s$, $D_t\mathbf{T}=0$. Differentiating with respect to $s$ and using the
preceding proposition,
\[
 0=D_sD_t\mathbf{T}=D_tD_s\mathbf{T}+\mathbf{R}(\mathbf{S},\mathbf{T})\mathbf{T}=D_t^2\mathbf{S}+\mathbf{R}(\mathbf{S},\mathbf{T})\mathbf{T}.
\]
Restriction to $s=0$ gives the stated equation. The sign agrees with
the convention
$\mathbf{R}(\mathbf{X},\mathbf{Y})\mathbf{Z}=\nabla_{\mathbf{X}}\nabla_{\mathbf{Y}}\mathbf{Z}-\nabla_{\mathbf{Y}}\nabla_{\mathbf{X}}\mathbf{Z}-\nabla_{[\mathbf{X},\mathbf{Y}]}\mathbf{Z}$
fixed in Definition~\ref{def:variedades-riemannianas-variedad-5}.
\end{proof}

\begin{proposition}[Existence and uniqueness]
\label{prop:existencia-unicidad-campos-jacobi}
Let \((M,\mathbf g)\) be a Riemannian manifold without boundary, let
\(\gamma\colon I\longrightarrow M\) be a geodesic, let \(t_0\in I\),
and let \(a,b\in T_{\gamma(t_0)}M\).
There is a unique Jacobi field \(\mathbf J\) along \(\gamma\) such that
\[
 \mathbf J(t_0)=a,\qquad \frac{D\mathbf J}{dt}(t_0)=b.
\]
The space of Jacobi fields along \(\gamma\) is a real vector space
of dimension \(2\dim M\).
\end{proposition}
\begin{proof}
Let \(n=\dim M\). Parallel transport of an orthonormal basis of
\(T_{\gamma(t_0)}M\) gives a parallel orthonormal frame
\((\mathbf E_1,\ldots,\mathbf E_n)\) along all of \(I\), by
Theorem~\ref{teo: existencia y unicidad transporte paralelo}.
Write
\[
 \mathbf J(t)=\displaystyle\sum_{j=1}^n u^j(t)\mathbf E_j(t),
 \qquad
 K^i_j(t):=\mathbf g_{\gamma(t)}
 \bigl(\mathbf R(\mathbf E_j,\dot\gamma)\dot\gamma,\mathbf E_i\bigr).
\]
Since \(D_t\mathbf E_j=0\), the Jacobi equation is equivalent to
\[
 (u^i)''(t)+\displaystyle\sum_{j=1}^n K^i_j(t)u^j(t)=0,
 \qquad 1\leq i\leq n.
\]
The functions \(K^i_j\) are smooth. Introducing
\(v^i=(u^i)'\) gives a first-order homogeneous linear system in
\(\mathbb R^{2n}\), with matrix
\[
 \begin{pmatrix}0&I_n\\-K(t)&0\end{pmatrix}.
\]
Theorem~\ref{teo: existencia y unicidad para ODEs lineales} gives a
unique solution on all of \(I\) for the components of the prescribed
data \(a,b\). Reconstructing the field in the parallel frame gives the
required solution. Linearity of the system shows that
\[
 \mathbf J\longmapsto
 \left(\mathbf J(t_0),\frac{D\mathbf J}{dt}(t_0)\right)
\]
is an isomorphism onto
\(T_{\gamma(t_0)}M\oplus T_{\gamma(t_0)}M\), proving the
dimension statement.
\end{proof}

\begin{proposition}[Representation by geodesic variations]
\label{prop:jacobi-representado-variacion-geodesica}
Let \((M,\mathbf g)\) be a Riemannian manifold without boundary, let
\(\gamma\colon I\longrightarrow M\) be a geodesic, and let
\(\mathbf J\) be a Jacobi field along it.
On every compact subinterval \(K\subseteq I\), the field
\(\mathbf J|_K\) is the variational field of a variation through
geodesics.
\end{proposition}
\begin{proof}
Fix \(t_0\in K\), set \(p=\gamma(t_0)\), and write
\[
 a=\mathbf J(t_0),\qquad
 b=\frac{D\mathbf J}{dt}(t_0),\qquad v=\dot\gamma(t_0).
\]
For small \(|s|\), let \(\sigma(s)=\exp_p(sa)\).
Denote by \(P_s\colon T_pM\longrightarrow T_{\sigma(s)}M\) parallel
transport along \(\sigma\) from \(0\) to \(s\), and define
\[
 \mathbf V(s):=P_s(v+sb).
\]
In a parallel frame along \(\sigma\), the components of
\(\mathbf V\) are those of \(v+sb\). Hence
\[
 \sigma(0)=p,\quad \sigma'(0)=a,\quad
 \mathbf V(0)=v,\quad D_s\mathbf V(0)=b.
\]
Consider
\[
 \alpha(s,t):=
 \exp_{\sigma(s)}\bigl((t-t_0)\mathbf V(s)\bigr).
\]
There is a single \(\delta>0\) for which this expression is defined on
\((-\delta,\delta)\times K\). Indeed, the domain of the geodesic
flow is open and contains the compact set
\[
 \{(t-t_0,(p,v))\mid t\in K\}.
\]
A finite cover by product neighborhoods yields a neighborhood of
\((p,v)\) common to all these times. For \(|s|<\delta\), the initial
data \((\sigma(s),\mathbf V(s))\) lie in that neighborhood.
Smoothness of the flow also proves that \(\alpha\) is smooth.

Each curve \(\alpha_s\) is geodesic and \(\alpha_0=\gamma|_K\).
If \(\widetilde{\mathbf J}\) is its variational field, then
\(\widetilde{\mathbf J}(t_0)=\sigma'(0)=a\).
By Lemma~\ref{lem:conmutacion-velocidades-variacion},
\[
 \frac{D\widetilde{\mathbf J}}{dt}(t_0)
 =\left.D_s\partial_t\alpha(s,t_0)\right|_{s=0}
 =D_s\mathbf V(0)=b.
\]
The Jacobi equation for geodesic variations and uniqueness for its
initial value problem give
\(\widetilde{\mathbf J}=\mathbf J|_K\).
\end{proof}

\begin{proposition}[Differential of the exponential map]
\label{prop:diferencial-exponencial-campo-jacobi}
\index{exponential map!differential}\index{Jacobi field!and exponential map}
Let \((M,\mathbf g)\) be a Riemannian manifold without boundary, let
\(p\in M\), let \(v\in\mathcal E_p\), and let \(w\in T_pM\).
Along \(\gamma_v(t)=\exp_p(tv)\), for \(0\leq t\leq1\), define
\[
 \mathbf J_w(t):=d(\exp_p)_{tv}(tw),
\]
using the identification \(T_{tv}(T_pM)\cong T_pM\).
Then \(\mathbf J_w\) is the Jacobi field with initial data
\[
 \mathbf J_w(0)=0,\qquad
 \frac{D\mathbf J_w}{dt}(0)=w,
\]
and, in particular, \(d(\exp_p)_v(w)=\mathbf J_w(1)\).
\end{proposition}
\begin{proof}
Since \(\mathcal E_p\) is open and star-shaped, for sufficiently small
\(|s|\) the family
\[
 \alpha(s,t)=\exp_p(t(v+sw)),\qquad 0\leq t\leq1,
\]
is a variation through geodesics. By the chain rule, its variational
field is \(\mathbf J_w(t)\).
All curves start at \(p\), so \(\mathbf J_w(0)=0\).
Moreover, \(\partial_t\alpha(s,0)=v+sw\) is a curve in the fixed
vector space \(T_pM\). Commutation of the velocities gives
\[
 \frac{D\mathbf J_w}{dt}(0)
 =\left.D_s\partial_t\alpha(s,0)\right|_{s=0}
 =\left.\frac{d}{ds}(v+sw)\right|_{s=0}=w.
\]
Theorem~\ref{teo:ecuacion-jacobi-variacion-geodesica} and
Proposition~\ref{prop:existencia-unicidad-campos-jacobi} identify
this field as the unique solution with those data. Evaluating its
definition at \(t=1\) gives the final identity.
\end{proof}

\begin{definition}[Conjugate points]
\label{def:puntos-conjugados-geodesica}
\index{conjugate point}\index{multiplicity of a conjugate point}
Let $(M,\mathbf{g})$ be a Riemannian manifold without boundary, let
$\gamma\colon I\to M$ be a geodesic, and let $t_0,t_1\in I$ be distinct. The points
$\gamma(t_0)$ and $\gamma(t_1)$ are \textbf{conjugate along
$\gamma$} if there exists a nonzero Jacobi field $\mathbf{J}$ such that
$\mathbf{J}(t_0)=\mathbf{J}(t_1)=0$. Their multiplicity is the dimension of the space of such
fields.
\end{definition}

\begin{proposition}[Conjugacy and critical values of the exponential map]
\label{prop:conjugados-valores-criticos-exponencial}
Let $(M,\mathbf{g})$ be a Riemannian manifold without boundary, let $p\in M$, and
let $v\ne0$ be a vector in the domain of $\exp_p$.
Then $\exp_p(v)$ is conjugate to $p$ along
$t\mapsto\exp_p(tv)$ if and only if $v$ is a critical point of $\exp_p$.
The multiplicity of the conjugate point is
$\dim\ker d(\exp_p)_v$.
\end{proposition}

\begin{proof}
By the preceding proposition, the map
$w\mapsto \mathbf{J}_w$ identifies $\ker d(\exp_p)_v$ with the Jacobi fields that
vanish at $t=0$ and $t=1$. Since the domain and codomain of
$d(\exp_p)_v$ have the same dimension, it is singular exactly
when its kernel is nontrivial.
\end{proof}

\begin{proposition}[Tangential and normal components of a Jacobi field]
\label{prop:jacobi-componentes-tangencial-normal}
Let \(\gamma\colon I\longrightarrow M\) be a unit-speed geodesic
and let \(\mathbf J\) be a Jacobi field. For \(t_0\in I\), set
\[
 a:=\mathbf g_{\gamma(t_0)}(\mathbf J(t_0),\dot\gamma(t_0)),
 \qquad
 b:=\mathbf g_{\gamma(t_0)}(D_t\mathbf J(t_0),\dot\gamma(t_0)).
\]
Then
\[
 \mathbf J^\top(t)=(a+b(t-t_0))\dot\gamma(t),
 \qquad \mathbf J^\perp:=\mathbf J-\mathbf J^\top
\]
are tangential and normal Jacobi fields, respectively.
In particular, a Jacobi field vanishing at two distinct times is normal.
\end{proposition}
\begin{proof}
Let \(h(t)=\mathbf g_{\gamma(t)}(\mathbf J(t),\dot\gamma(t))\).
Since \(D_t\dot\gamma=0\), the Jacobi equation and the curvature
symmetries imply
\[
 h''(t)
 =\mathbf g_{\gamma(t)}(D_t^2\mathbf J,\dot\gamma)
 =-\mathbf g_{\gamma(t)}
   (\mathbf R(\mathbf J,\dot\gamma)\dot\gamma,\dot\gamma)=0.
\]
The initial values \(h(t_0)=a\), \(h'(t_0)=b\) give
\(h(t)=a+b(t-t_0)\). Because
\(\mathbf R(\dot\gamma,\dot\gamma)\dot\gamma=0\),
the field \(h\dot\gamma\) satisfies the Jacobi equation. Its difference
from \(\mathbf J\) also satisfies that equation and is orthogonal
to \(\dot\gamma\).
If \(\mathbf J\) vanishes at two distinct times, the affine function
\(h\) has two zeros and therefore vanishes identically.
\end{proof}

\begin{proposition}[The two-point boundary problem for Jacobi fields]
\label{prop:jacobi-valores-dos-extremos}
Let \(\gamma\colon[a,b]\longrightarrow M\) be a geodesic, with
\(a<b\). If \(\gamma(a)\) and \(\gamma(b)\) are not conjugate along
it, then for every \(v\in T_{\gamma(a)}M\) and
\(w\in T_{\gamma(b)}M\) there is a unique Jacobi field
\(\mathbf J\) such that \(\mathbf J(a)=v\) and \(\mathbf J(b)=w\).
\end{proposition}
\begin{proof}
For \(u\in T_{\gamma(a)}M\), let \(\mathbf J_u\) be the Jacobi
field with \(\mathbf J_u(a)=0\) and \(D_t\mathbf J_u(a)=u\).
The linear map
\[
 B\colon T_{\gamma(a)}M\longrightarrow T_{\gamma(b)}M,
 \qquad B(u):=\mathbf J_u(b),
\]
is injective: a nonzero element of its kernel would give a nonzero
Jacobi field vanishing at both endpoints. Since the dimensions agree,
\(B\) is an isomorphism.
Let \(\mathbf H\) be the Jacobi field with initial data
\(\mathbf H(a)=v\), \(D_t\mathbf H(a)=0\). Then
\[
 \mathbf J:=\mathbf H+\mathbf J_{B^{-1}(w-\mathbf H(b))}
\]
has the prescribed endpoint values. The difference of two solutions
would vanish at both endpoints and hence would be zero.
\end{proof}

To formulate the constant-curvature models, define
\[
s_c(t):=
\begin{cases}
\dfrac{\sin(\sqrt c\,t)}{\sqrt c},&c>0,\\[6pt]
t,&c=0,\\[4pt]
\dfrac{\sinh(\sqrt{-c}\,t)}{\sqrt{-c}},&c<0,
\end{cases}
\qquad c_c(t):=s_c'(t).
\]
These functions satisfy $s_c''+cs_c=0$,
$s_c(0)=0$, and $s_c'(0)=1$.

\begin{proposition}[Jacobi fields in constant curvature]
\label{prop:jacobi-curvatura-seccional-constante}
Let $(M,\mathbf{g})$ be a Riemannian manifold without boundary of constant
sectional curvature $c$, let $\gamma$ be a unit-speed geodesic, and let
$\mathbf{J}$ be a Jacobi field along $\gamma$. Let $P_t$ be parallel transport
from $\gamma(0)$ to $\gamma(t)$. If
$a=\mathbf{J}(0)$ and $b=D_t\mathbf{J}(0)$, decompose
$a=a^\top+a^\perp$ and $b=b^\top+b^\perp$ with respect to
$\mathbb R\dot\gamma(0)\oplus\dot\gamma(0)^\perp$. Then
\[
 \mathbf{J}(t)=P_t\bigl(a^\top+t b^\top+c_c(t)a^\perp+s_c(t)b^\perp\bigr).
\]
In particular, if $\mathbf{J}(0)=0$ and $D_t\mathbf{J}(0)=w\perp\dot\gamma(0)$, then
$\mathbf{J}(t)=s_c(t)P_tw$.
\end{proposition}
\begin{proof}
In constant curvature,
$\mathbf{R}(\mathbf{X},\dot\gamma)\dot\gamma=c\mathbf{X}$ for
$\mathbf{X}\perp\dot\gamma$, while the curvature term vanishes in the
direction of $\dot\gamma$. In a parallel frame, the Jacobi equation thus
splits into $u''=0$ for the tangential component and
$u''+cu=0$ for each normal component. The solutions with the given initial
data are exactly those in the formula.
\end{proof}

To formulate the variational argument underlying comparison,
we introduce the index form. A partition
$a=t_0<t_1<\cdots<t_N=b$ is \emph{admissible} for a field $\mathbf{V}$ along
$\gamma$ if $\mathbf{V}$ is continuous on $[a,b]$ and its restriction to each
$[t_{r-1},t_r]$ is smooth up to the endpoints. Denote
the space of these piecewise smooth fields by $\mathcal V_{\mathrm{ps}}(\gamma)$
and write
\[
\mathcal V_0(\gamma)
:=\{\mathbf{V}\in\mathcal V_{\mathrm{ps}}(\gamma)\mid \mathbf{V}(a)=\mathbf{V}(b)=0\}.
\]

\begin{definition}[Index form]\label{def:forma-indice}\index{index form}
Let $(M,\mathbf{g})$ be a Riemannian manifold without boundary and let
$\gamma\colon[a,b]\longrightarrow M$ be a geodesic. For
$\mathbf{V},\mathbf{W}\in\mathcal V_{\mathrm{ps}}(\gamma)$, define
\[
I_\gamma(\mathbf{V},\mathbf{W})
:=
\int_a^b
\left[
\left\langle\frac{D\mathbf{V}}{dt},\frac{D\mathbf{W}}{dt}\right\rangle_{\mathbf{g}}
-\left\langle \mathbf{R}(\mathbf{V},\dot\gamma)\dot\gamma,\mathbf{W}\right\rangle_{\mathbf{g}}
\right]dt.
\]
Covariant derivatives are taken on each subinterval of an admissible
partition.
\end{definition}

\begin{proposition}[Properties and integration by parts for the index form]
\label{prop:integracion-por-partes-forma-indice}
\index{index form!integration by parts}
Let $(M,\mathbf{g})$ be a Riemannian manifold without boundary and let
$\gamma\colon[a,b]\to M$ be a geodesic. The form $I_\gamma$ is well-defined,
bilinear, and symmetric. If
$a=t_0<\cdots<t_N=b$ is a common admissible partition for $\mathbf{V}$ and $\mathbf{W}$, and
\[
\left[\frac{D\mathbf{V}}{dt}\right]_{t_r}
:=
\frac{D\mathbf{V}}{dt}(t_r^+)-\frac{D\mathbf{V}}{dt}(t_r^-),
\qquad 1\leq r\leq N-1,
\]
then
\begin{align}
I_\gamma(\mathbf{V},\mathbf{W})
={}&
\left[
\left\langle\frac{D\mathbf{V}}{dt},\mathbf{W}\right\rangle_{\mathbf{g}}
\right]_{a}^{b}
-\displaystyle\sum_{r=1}^{N-1}
\left\langle
\left[\frac{D\mathbf{V}}{dt}\right]_{t_r},\mathbf{W}(t_r)
\right\rangle_{\mathbf{g}}
\label{eq:integracion-partes-forma-indice}\\
&-
\int_a^b
\left\langle
\frac{D^2\mathbf{V}}{dt^2}+\mathbf{R}(\mathbf{V},\dot\gamma)\dot\gamma,
\mathbf{W}
\right\rangle_{\mathbf{g}}dt.\notag
\end{align}
The values at $a$ and $b$ are the limits from the interior. In particular,
if $\mathbf{J}$ is a smooth Jacobi field, then
\begin{equation}
\label{eq:forma-indice-jacobi-borde}
I_\gamma(\mathbf{J},\mathbf{W})
=
\left[
\left\langle\frac{D\mathbf{J}}{dt},\mathbf{W}\right\rangle_{\mathbf{g}}
\right]_{a}^{b}.
\end{equation}
\end{proposition}

\begin{proof}
On a common admissible partition, the integrand is continuous in the interior
of each subinterval and has one-sided limits at its endpoints. Refining the
partition merely subdivides existing integrals, and the partition points
have measure zero; the value is therefore independent of the partition.
Bilinearity is immediate. The symmetries of the curvature tensor imply
\[
\langle \mathbf{R}(\mathbf{V},\dot\gamma)\dot\gamma,\mathbf{W}\rangle_{\mathbf{g}}
=\langle \mathbf{R}(\mathbf{W},\dot\gamma)\dot\gamma,\mathbf{V}\rangle_{\mathbf{g}},
\]
so $I_\gamma$ is symmetric.

On each interval $[t_{r-1},t_r]$, we have
\[
\frac{d}{dt}
\left\langle\frac{D\mathbf{V}}{dt},\mathbf{W}\right\rangle_{\mathbf{g}}
=
\left\langle\frac{D^2\mathbf{V}}{dt^2},\mathbf{W}\right\rangle_{\mathbf{g}}
+
\left\langle\frac{D\mathbf{V}}{dt},\frac{D\mathbf{W}}{dt}\right\rangle_{\mathbf{g}}.
\]
Integrate this equality and sum over $r$. At an interior point $t_r$,
the boundary terms contribute
\[
\left\langle\frac{D\mathbf{V}}{dt}(t_r^-),\mathbf{W}(t_r)\right\rangle_{\mathbf{g}}
-
\left\langle\frac{D\mathbf{V}}{dt}(t_r^+),\mathbf{W}(t_r)\right\rangle_{\mathbf{g}}
=-
\left\langle
\left[\frac{D\mathbf{V}}{dt}\right]_{t_r},\mathbf{W}(t_r)
\right\rangle_{\mathbf{g}},
\]
because $\mathbf{W}$ is continuous. Including the curvature term gives
\eqref{eq:integracion-partes-forma-indice}. If $\mathbf{J}$ is a Jacobi field,
it has no jumps and the integrand in the second line vanishes, giving
\eqref{eq:forma-indice-jacobi-borde}.
\end{proof}

The index form is the second variation of energy at a geodesic with
fixed endpoints.
The geodesic equation also has a variational interpretation.
For a smooth curve \(c\colon[a,b]\longrightarrow M\), with \(a<b\),
define its \textbf{energy} by
\[
 E(c):=\frac12\int_a^b|\dot c(t)|_{\mathbf g_{c(t)}}^2\,dt.
\]
Unlike length, energy depends on the parametrization.
The Cauchy--Schwarz inequality gives
\[
 L_{\mathbf g}(c)^2
 \leq(b-a)\int_a^b|\dot c(t)|_{\mathbf g_{c(t)}}^2\,dt
 =2(b-a)E(c),
\]
with equality if and only if the speed is constant.
The same formulas apply to continuous, piecewise smooth curves.

\begin{proposition}[First variation of energy]
\label{prop:primera-variacion-energia}
Let \(\alpha\colon(-\delta,\delta)\times[a,b]\longrightarrow M\)
be a smooth variation, and set
\(\mathbf S=d\alpha(\boldsymbol{\partial}_s)\),
\(\mathbf T=d\alpha(\boldsymbol{\partial}_t)\). Then
\[
 \frac{d}{ds}E(\alpha_s)
 =\left[\mathbf g_\alpha(\mathbf S,\mathbf T)\right]_{t=a}^{t=b}
 -\int_a^b\mathbf g_{\alpha(s,t)}(\mathbf S,D_t\mathbf T)\,dt.
\]
A smooth curve is geodesic if and only if its first variation of
energy vanishes for every smooth variation with fixed endpoints.
\end{proposition}
\begin{proof}
Smoothness and compactness of \([a,b]\) allow differentiation under
the integral sign. Metric compatibility and commutation of the
velocities give
\[
 \frac{d}{ds}E(\alpha_s)
 =\int_a^b\mathbf g_\alpha(D_s\mathbf T,\mathbf T)\,dt
 =\int_a^b\mathbf g_\alpha(D_t\mathbf S,\mathbf T)\,dt.
\]
Integrating
\(\partial_t\mathbf g_\alpha(\mathbf S,\mathbf T)
=\mathbf g_\alpha(D_t\mathbf S,\mathbf T)
+\mathbf g_\alpha(\mathbf S,D_t\mathbf T)\)
over \([a,b]\) proves the formula.
If \(\gamma=\alpha_0\) is geodesic and the endpoints are fixed, then
\(\mathbf S(s,a)=\mathbf S(s,b)=0\) and \(D_t\dot\gamma=0\);
the derivative at \(s=0\) therefore vanishes.

Conversely, every smooth field \(\mathbf V\) along \(\gamma\)
vanishing at the endpoints is realized by
\[
 \alpha(s,t):=\exp_{\gamma(t)}(s\mathbf V(t)).
\]
This expression is defined for a uniform interval of values of
\(s\): the domain of the exponential map is open and contains the
zero section over the compact set \(\gamma([a,b])\).
Its variational field is \(\mathbf V\), since
\(d(\exp_{\gamma(t)})_0\) is the identity.
If \(D_t\dot\gamma\) were nonzero at an interior time, choose a
smooth function \(\chi\geq0\), compactly supported in \((a,b)\)
and positive near that time, and take
\(\mathbf V=\chi D_t\dot\gamma\). The first variation would be
\[
 -\int_a^b\chi(t)|D_t\dot\gamma(t)|_{\mathbf g_{\gamma(t)}}^2\,dt<0,
\]
a contradiction. Thus \(D_t\dot\gamma=0\) on \((a,b)\), and also
at the endpoints by continuity.
\end{proof}

\begin{proposition}[Second variation of energy]
\label{prop:segunda-variacion-energia-forma-indice}
\index{second variation!of energy}
Let $(M,\mathbf{g})$ be a Riemannian manifold without boundary and let
$\alpha\colon(-\varepsilon,\varepsilon)\times[a,b]\to M$ be a
smooth variation of a geodesic $\gamma=\alpha(0,\cdot)$ with fixed
endpoints, and let $\mathbf{V}$ be its variational field. If
\[
E(\alpha_s):=\frac12\int_a^b
\bigl|d\alpha_{(s,t)}(\boldsymbol{\partial}_t)\bigr|_{\mathbf{g}}^2dt,
\]
then
\[
\left.\frac{d^2}{ds^2}\right|_{s=0}E(\alpha_s)=I_\gamma(\mathbf{V},\mathbf{V}).
\]
Without the fixed-endpoint hypothesis, the formula also contains the boundary
term
\[
 \left.
 \left[\left\langle D_s\mathbf{S},\mathbf{T}\right\rangle_{\mathbf{g}}\right]_a^b
 \right|_{s=0},
 \qquad
 \mathbf{S}=d\alpha(\boldsymbol{\partial}_s),\quad \mathbf{T}=d\alpha(\boldsymbol{\partial}_t).
\]
\end{proposition}
\begin{proof}
Write $\mathbf{S}=d\alpha(\boldsymbol{\partial}_s)$ and $\mathbf{T}=d\alpha(\boldsymbol{\partial}_t)$.
Metric compatibility and the identity $D_s\mathbf{T}=D_t\mathbf{S}$ give
\[
\frac{d}{ds}E(\alpha_s)=\int_a^b\langle D_t\mathbf{S},\mathbf{T}\rangle_{\mathbf{g}}dt.
\]
Differentiating again and restricting to $s=0$, use
$D_sD_t\mathbf{S}=D_tD_s\mathbf{S}+\mathbf{R}(\mathbf{S},\mathbf{T})\mathbf{S}$ and $D_s\mathbf{T}=D_t\mathbf{S}$ to obtain
\[
\left.\frac{d^2}{ds^2}\right|_{s=0}E(\alpha_s)
=\int_a^b
\left[
\langle D_tD_s\mathbf{S},\dot\gamma\rangle_{\mathbf{g}}
+\langle \mathbf{R}(\mathbf{V},\dot\gamma)\mathbf{V},\dot\gamma\rangle_{\mathbf{g}}
+\left|\frac{D\mathbf{V}}{dt}\right|_{\mathbf{g}}^2
\right]dt.
\]
After integration by parts, the first integral reduces to the boundary term
in the statement, since the geodesic satisfies
$D_t\dot\gamma=0$. Because the endpoints of the variation are fixed,
$D_s\mathbf{S}=0$ there, and this term vanishes. Finally,
\[
\langle \mathbf{R}(\mathbf{V},\dot\gamma)\mathbf{V},\dot\gamma\rangle_{\mathbf{g}}
=-\langle \mathbf{R}(\mathbf{V},\dot\gamma)\dot\gamma,\mathbf{V}\rangle_{\mathbf{g}},
\]
which gives the asserted formula.
\end{proof}

The second variation formula remains valid for a continuous variation
that is smooth on each rectangle
\((-\delta,\delta)\times[t_{j-1},t_j]\) of a fixed partition.
Indeed, differentiate and integrate by parts on each rectangle.
The values of \(\mathbf S\) and \(D_s\mathbf S\) agree on common
sides, since they are obtained by differentiating the same curve
\(s\mapsto\alpha(s,t_j)\) with respect to \(s\).
At a smooth geodesic \(\alpha_0\), the interior boundary terms
cancel. The result is the sum of the integrals defining \(I_\gamma\).

\begin{proposition}[Nonnegativity of the index form along a minimizer]
\label{prop:indice-no-negativo-minimizante}
If \(\gamma\colon[a,b]\longrightarrow M\) is a unit-speed
minimizing geodesic, with \(a<b\), then
\[
 I_\gamma(\mathbf V,\mathbf V)\geq0
 \qquad\text{for every }\mathbf V\in\mathcal V_0(\gamma).
\]
\end{proposition}
\begin{proof}
For an admissible partition of \(\mathbf V\), the expression
\(\alpha(s,t)=\exp_{\gamma(t)}(s\mathbf V(t))\) defines, for
sufficiently small \(|s|\), a continuous variation that is smooth
on each rectangle and has fixed endpoints.
Compactness and finiteness of the partition give an interval of
values of \(s\) common to all rectangles.
Since \(|\partial_t\alpha(0,t)|_{\mathbf g}=1\), after shrinking
this interval we can require
\(|\partial_t\alpha(s,t)|_{\mathbf g}>1/2\) on each rectangle.
Thus \(\alpha_s\) is admissible.

The minimizing property and the inequality between length and
energy give
\[
 E(\alpha_s)\geq
 \frac{L_{\mathbf g}(\alpha_s)^2}{2(b-a)}
 \geq\frac{d_{\mathbf g}(\gamma(a),\gamma(b))^2}{2(b-a)}
 =\frac{b-a}{2}=E(\gamma).
\]
The function \(s\mapsto E(\alpha_s)\) has a minimum at zero.
By the second variation formula, its second derivative there is
\(I_\gamma(\mathbf V,\mathbf V)\), which must be nonnegative.
\end{proof}

\begin{theorem}[A geodesic does not minimize past a conjugate point]
\label{teo:conjugado-interior-impide-minimizacion}
Let \(\gamma\colon[a,b]\longrightarrow M\) be a unit-speed geodesic.
If \(\gamma(a)\) is conjugate to \(\gamma(t_0)\) along \(\gamma\)
for some \(a<t_0<b\), there is
\(\mathbf X\in\mathcal V_0(\gamma)\) such that
\(I_\gamma(\mathbf X,\mathbf X)<0\).
In particular, \(\gamma\) is not minimizing.
\end{theorem}
\begin{proof}
Let \(\mathbf J\) be a nonzero Jacobi field on \([a,t_0]\) with
\(\mathbf J(a)=\mathbf J(t_0)=0\). Define
\[
 \mathbf V(t):=
 \begin{cases}
 \mathbf J(t),&a\leq t\leq t_0,\\
 0,&t_0\leq t\leq b.
 \end{cases}
\]
This field is continuous, piecewise smooth, and zero at the endpoints.
Set \(u=D_t\mathbf J(t_0)\). Uniqueness for the Jacobi equation
implies \(u\neq0\), since otherwise both initial data of
\(\mathbf J\) at \(t_0\) would be zero.
Choose \(\chi\in C_c^\infty((a,b))\) with \(\chi(t_0)=1\), and
define
\[
 \mathbf W(t):=\chi(t)P_{t_0,t}u,
\]
where \(P_{t_0,t}\) denotes parallel transport along \(\gamma\).
Integration by parts for the index form, or its direct application
on \([a,t_0]\) and \([t_0,b]\), gives
\[
 I_\gamma(\mathbf V,\mathbf V)=0,\qquad
 I_\gamma(\mathbf V,\mathbf W)
 =\mathbf g_{\gamma(t_0)}(u,\mathbf W(t_0))=|u|_{\mathbf g_{\gamma(t_0)}}^2.
\]
The sign of the second expression corresponds to the jump
\([D_t\mathbf V]_{t_0}=-u\).
For \(\eta>0\), let \(\mathbf X_\eta=\mathbf V-\eta\mathbf W\).
Then
\[
 I_\gamma(\mathbf X_\eta,\mathbf X_\eta)
 =-2\eta|u|_{\mathbf g_{\gamma(t_0)}}^2
   +\eta^2I_\gamma(\mathbf W,\mathbf W)<0
\]
for sufficiently small \(\eta\).
The preceding proposition rules out minimization by \(\gamma\).
\end{proof}

\begin{lemma}[Index lemma]\label{lem:indice-campos-jacobi}
\index{index lemma}
Let $(M,\mathbf{g})$ be a Riemannian manifold without boundary, let
$\gamma\colon[a,b]\to M$ be a geodesic, and suppose that
$\gamma(a)$ is not conjugate to $\gamma(t)$ along $\gamma$ for
any $t\in(a,b]$. Let $\mathbf{J}$ be a Jacobi field along $\gamma$
with $\mathbf{J}(a)=0$. If $\mathbf{V}\in\mathcal V_{\mathrm{ps}}(\gamma)$ satisfies
\[
\mathbf{V}(a)=0,
\qquad
\mathbf{V}(b)=\mathbf{J}(b),
\]
then
\[
I_\gamma(\mathbf{J},\mathbf{J})\leq I_\gamma(\mathbf{V},\mathbf{V}).
\]
Equality holds if and only if $\mathbf{V}=\mathbf{J}$.
\end{lemma}

\begin{proof}
Choose an orthonormal basis $(\mathbf{e}_1,\dots,\mathbf{e}_n)$ of
$T_{\gamma(a)}M$ and extend it to a parallel frame along
$\gamma$. For each $i$, let $\mathbf{J}_i$ be the unique Jacobi field determined
by
\[
\mathbf{J}_i(a)=0,
\qquad
\frac{D\mathbf{J}_i}{dt}(a)=\mathbf{e}_i.
\]
For $t>a$, the vectors $\mathbf{J}_1(t),\dots,\mathbf{J}_n(t)$ form a basis of
$T_{\gamma(t)}M$. Indeed, a linear combination vanishing at
$t$ would give a Jacobi field vanishing at $a$ and $t$; by the
absence of conjugate points, it would have to be the zero field, and its initial
data would force all the coefficients to vanish.

We first prove that
\begin{equation}
\label{eq:positividad-forma-indice-campos-propios}
I_\gamma(\mathbf{W},\mathbf{W})\geq0
\qquad\text{for every }\mathbf{W}\in\mathcal V_0(\gamma),
\end{equation}
and determine the equality case. For $t>a$, write uniquely
\[
\mathbf{W}(t)=f^i(t)\mathbf{J}_i(t).
\]
The functions $f^i$ are piecewise smooth. They also extend smoothly
from the right to $a$ on the first admissible subinterval. To see this,
express the fields in the parallel frame and denote by $A(t)$ the
matrix whose columns are the components of the $\mathbf{J}_i(t)$. Since
$A(a)=0$ and $A'(a)=I$, the integral formula
\[
A(t)=(t-a)B(t),
\qquad
B(t):=\int_0^1A'\bigl(a+s(t-a)\bigr)ds,
\]
defines a smooth matrix with $B(a)=I$. Similarly, $\mathbf{W}(a)=0$ allows us to
write $\mathbf{W}(t)=(t-a)w(t)$ with $w$ smooth on the first subinterval. For
$t$ close to $a$, $B(t)$ is invertible and
\[
(f^1(t),\dots,f^n(t))^{\mathsf T}=B(t)^{-1}w(t),
\]
which proves the extension and controls the endpoint where all the
$\mathbf{J}_i$ vanish.

For each pair $i,j$, the function
\[
\left\langle\frac{D\mathbf{J}_i}{dt},\mathbf{J}_j\right\rangle_{\mathbf{g}}
-\left\langle \mathbf{J}_i,\frac{D\mathbf{J}_j}{dt}\right\rangle_{\mathbf{g}}
\]
is constant. Its derivative vanishes by the Jacobi equation and symmetry
of the endomorphism
$\mathbf{Z}\mapsto \mathbf{R}(\mathbf{Z},\dot\gamma)\dot\gamma$; its value at $a$ is zero. Therefore,
\begin{equation}
\label{eq:simetria-wronskiano-jacobi-indice}
\left\langle\frac{D\mathbf{J}_i}{dt},\mathbf{J}_j\right\rangle_{\mathbf{g}}
=\left\langle \mathbf{J}_i,\frac{D\mathbf{J}_j}{dt}\right\rangle_{\mathbf{g}}.
\end{equation}

On each subinterval where $\mathbf{W}$ is smooth, set
\[
\mathbf{X}:=f^i\frac{D\mathbf{J}_i}{dt},
\qquad
\mathbf{Y}:=\frac{df^i}{dt}\mathbf{J}_i.
\]
Then $\frac{D\mathbf{W}}{dt}=\mathbf{X}+\mathbf{Y}$. The Jacobi equation and
\eqref{eq:simetria-wronskiano-jacobi-indice} give
\[
\frac{D\mathbf{X}}{dt}
=\frac{df^i}{dt}\frac{D\mathbf{J}_i}{dt}
-\mathbf{R}(\mathbf{W},\dot\gamma)\dot\gamma,
\qquad
\left\langle \mathbf{W},
\frac{df^i}{dt}\frac{D\mathbf{J}_i}{dt}\right\rangle_{\mathbf{g}}
=\langle \mathbf{Y},\mathbf{X}\rangle_{\mathbf{g}}.
\]
Consequently,
\begin{align}
\left|\frac{D\mathbf{W}}{dt}\right|_{\mathbf{g}}^2
-\langle \mathbf{R}(\mathbf{W},\dot\gamma)\dot\gamma,\mathbf{W}\rangle_{\mathbf{g}}
=
\frac{d}{dt}\langle \mathbf{W},\mathbf{X}\rangle_{\mathbf{g}}+|\mathbf{Y}|_{\mathbf{g}}^2.
\label{eq:identidad-fundamental-lema-indice}
\end{align}

The coefficients $f^i$ are continuous at the partition points because
$\mathbf{W}$ is continuous and $(\mathbf{J}_i(t))$ is a basis; therefore, $\mathbf{X}$ is also
continuous there. Integrating
\eqref{eq:identidad-fundamental-lema-indice} over all the
subintervals, the interior boundary terms cancel. At $a$, they
vanish because $\mathbf{W}(a)=0$ and $\mathbf{X}$ has a limit, and at $b$ they vanish because
$\mathbf{W}(b)=0$. Thus,
\[
I_\gamma(\mathbf{W},\mathbf{W})=\int_a^b|\mathbf{Y}|_{\mathbf{g}}^2dt\geq0.
\]
If equality holds, $\mathbf{Y}=0$ on each subinterval. Since the $\mathbf{J}_i(t)$ are
linearly independent for $t>a$, all the functions $f^i$ are
constant on each subinterval; their continuity shows that they are constant
on all of $[a,b]$. Then $\mathbf{W}$ is a Jacobi field. Since
$\mathbf{W}(a)=\mathbf{W}(b)=0$ and the endpoints are not conjugate, it follows that $\mathbf{W}=0$.
Conversely, $\mathbf{W}=0$ gives equality. This proves
\eqref{eq:positividad-forma-indice-campos-propios} and its equality case.

Return to the fields in the statement and set $\mathbf{W}:=\mathbf{V}-\mathbf{J}$. Then
$\mathbf{W}\in\mathcal V_0(\gamma)$. The integration-by-parts formula and the
Jacobi equation imply
\[
I_\gamma(\mathbf{J},\mathbf{W})
=\left[
\left\langle\frac{D\mathbf{J}}{dt},\mathbf{W}\right\rangle_{\mathbf{g}}
\right]_a^b=0.
\]
By symmetry and bilinearity,
\[
I_\gamma(\mathbf{V},\mathbf{V})
=I_\gamma(\mathbf{J},\mathbf{J})+2I_\gamma(\mathbf{J},\mathbf{W})+I_\gamma(\mathbf{W},\mathbf{W})
=I_\gamma(\mathbf{J},\mathbf{J})+I_\gamma(\mathbf{W},\mathbf{W}).
\]
The inequality already proved gives the conclusion, and equality holds
exactly when $\mathbf{W}=0$, that is, when $\mathbf{V}=\mathbf{J}$.
\end{proof}

\begin{theorem}[Comparison of Jacobi fields]
\label{teo: de comparacion de Rauch}
\index{comparison of Jacobi fields}
Let $(M,\mathbf{g})$ be a Riemannian manifold without boundary of dimension at
least two, let
$\gamma\colon[0,b]\to M$ be a unit-speed geodesic, and let $\mathbf{J}$ be a normal
Jacobi field, that is, $\langle \mathbf{J},\dot\gamma\rangle_{\mathbf{g}}=0$, with
$\mathbf{J}(0)=0$.
Define
\[
\rho_c:=
\begin{cases}
+\infty,&c\leq0,\\[2pt]
\dfrac{\pi}{\sqrt c},&c>0,
\end{cases}
\]
and let $\tau$ be the first time conjugate to $0$ along $\gamma$,
with the convention $\tau=+\infty$ if there is none.
\begin{enumerate}[label=(\alph*)]
\item If every sectional curvature of $M$ is less than or equal to $c$, then
\[
 |\mathbf{J}(t)|_{\mathbf{g}}\geq s_c(t)\left|\frac{D\mathbf{J}}{dt}(0)\right|_{\mathbf{g}}.
\]
The inequality holds for $0\leq t\leq\min\{b,\rho_c\}$. Moreover,
$\gamma(0)$ has no conjugate points along $\gamma$ as long as
$0<t<\min\{b,\rho_c\}$.
\item If every sectional curvature of $M$ is greater than or equal to $c$, then
\[
 |\mathbf{J}(t)|_{\mathbf{g}}\leq s_c(t)\left|\frac{D\mathbf{J}}{dt}(0)\right|_{\mathbf{g}}
\]
for $0\leq t\leq\min\{b,\tau,\rho_c\}$. If $c>0$ and
$b\geq\rho_c$, then $\tau\leq\rho_c$.
\end{enumerate}

If $\mathbf{J}$ is not the zero field and $0<t_0<\min\{b,\tau,\rho_c\}$ in the case
involving $\tau$, or $0<t_0<\min\{b,\rho_c\}$ in the other case,
equality in the corresponding comparison at $t_0$ holds if and
only if
\[
\mathbf{J}(t)=s_c(t)P_t\left(\frac{D\mathbf{J}}{dt}(0)\right)
\quad\text{and}\quad
\mathbf{R}(P_tw,\dot\gamma(t))\dot\gamma(t)=cP_tw
\]
for $0<t\leq t_0$, where
$w=\frac{D\mathbf{J}}{dt}(0)$ and $P_t$ is parallel transport from
$\gamma(0)$ to $\gamma(t)$.
\end{theorem}
\begin{proof}
Set $w=\frac{D\mathbf{J}}{dt}(0)$. If $w=0$, uniqueness for the Jacobi equation
implies $\mathbf{J}=0$, and both inequalities are immediate. Thus,
suppose that $w\ne0$. Choose a parallel orthonormal frame
$(\mathbf{E}_2,\dots,\mathbf{E}_n)$ of the orthogonal complement of $\dot\gamma$ and write
\[
\mathbf{J}(t)=u^\alpha(t)\mathbf{E}_\alpha(t),
\qquad
u(t)=(u^2(t),\dots,u^n(t))\in\mathbb R^{n-1}.
\]
In particular,
$|\mathbf{J}|_{\mathbf{g}}=|u|$, $\left|\frac{D\mathbf{J}}{dt}\right|_{\mathbf{g}}=|u'|$, and
$u'(0)$ represents $w$.

For $t<\rho_c$, consider on $[0,t]$ the model index form
\[
I_c(v,v):=\int_0^t\bigl(|v'(s)|^2-c|v(s)|^2\bigr)ds
\]
on fields with values in $\mathbb R^{n-1}$. The model Jacobi equation
is $v''+cv=0$. Since $s_c(s)>0$ for $0<s<\rho_c$, this equation
has no nonzero solution vanishing at $0$ and $t$. The
field
\[
\overline u_t(s):=\frac{s_c(s)}{s_c(t)}u(t)
\]
minimizes $I_c$ among fields taking the value $0$ at $s=0$ and $u(t)$ at
$s=t$. Indeed, for such a field, write
$v(s)=s_c(s)f(s)$ on $(0,t]$. The factorization at the initial endpoint
used in the proof of Lemma~\ref{lem:indice-campos-jacobi} shows that
$f$ extends regularly to $0$. Since $s_c''=-cs_c$,
\[
 |v'|^2-c|v|^2
 =\frac{d}{ds}\bigl(s_c c_c|f|^2\bigr)+s_c^2|f'|^2.
\]
The boundary term at $0$ vanishes. The value at $t$ is fixed and the
remaining integral is nonnegative, with equality exactly when $f$ is
constant. This proves the minimizing property and, for
$v=\overline u_t$, gives
\begin{equation}
\label{eq:indice-modelo-comparacion-jacobi}
I_c(\overline u_t,\overline u_t)
=\frac{c_c(t)}{s_c(t)}|u(t)|^2.
\end{equation}

First suppose that $K\leq c$. Since $\mathbf{J}$ is normal and $\gamma$ has
unit speed,
\[
\langle \mathbf{R}(\mathbf{J},\dot\gamma)\dot\gamma,\mathbf{J}\rangle_{\mathbf{g}}
\leq c|\mathbf{J}|_{\mathbf{g}}^2.
\]
The index lemma for the model and this inequality imply
\begin{align*}
\frac{c_c(t)}{s_c(t)}|\mathbf{J}(t)|_{\mathbf{g}}^2
&=I_c(\overline u_t,\overline u_t)\\
&\leq I_c(u,u)\\
&\leq I_{\gamma|_{[0,t]}}(\mathbf{J},\mathbf{J})\\
&=\left\langle\frac{D\mathbf{J}}{dt}(t),\mathbf{J}(t)\right\rangle_{\mathbf{g}}.
\end{align*}
The last equality is
\eqref{eq:forma-indice-jacobi-borde}, since $\mathbf{J}(0)=0$. On every interval
where $\mathbf{J}(t)\ne0$, it follows that
\[
\frac{d}{dt}\log\left(\frac{|\mathbf{J}(t)|_{\mathbf{g}}}{s_c(t)}\right)
=
\frac{\left\langle\frac{D\mathbf{J}}{dt},\mathbf{J}\right\rangle_{\mathbf{g}}}{|\mathbf{J}|_{\mathbf{g}}^2}
-\frac{c_c}{s_c}
\geq0.
\]
In the parallel frame, $u(0)=0$ and $u'(0)=w$; hence,
\[
\lim_{t\to 0^{+}}\frac{|\mathbf{J}(t)|_{\mathbf{g}}}{s_c(t)}
=\lim_{t\to 0^{+}}
\frac{\left|\frac{u(t)}{t}\right|}{\frac{s_c(t)}{t}}=|w|_{\mathbf{g}}.
\]
The quotient is nondecreasing, giving
\[
|\mathbf{J}(t)|_{\mathbf{g}}\geq s_c(t)|w|_{\mathbf{g}}
\]
as long as $s_c(t)>0$. In particular, $\mathbf{J}$ cannot vanish there. If an
arbitrary Jacobi field vanishes at $0$ and at some $t>0$, its tangential
component has the form $at\dot\gamma$ and must vanish; the field is
therefore normal. The preceding inequality excludes such a nonzero field for
$t<\rho_c$. The comparison extends to the endpoint
$t=\min\{b,\rho_c\}$ by continuity. This proves (a).

Now suppose that $K\geq c$ and fix
$0<t<\min\{b,\tau,\rho_c\}$. Along $\gamma|_{[0,t]}$, define
the field
\[
\mathbf{V}_t(s):=\frac{s_c(s)}{s_c(t)}P_{t,s}\mathbf{J}(t),
\]
where $P_{t,s}$ is parallel transport from $\gamma(t)$ to
$\gamma(s)$. The fields $\mathbf{V}_t$ and $\mathbf{J}$ have the same values at the
endpoints. Since there are no conjugate points before $t$,
Lemma~\ref{lem:indice-campos-jacobi} and the lower curvature bound give
\begin{align*}
\left\langle\frac{D\mathbf{J}}{dt}(t),\mathbf{J}(t)\right\rangle_{\mathbf{g}}
&=I_{\gamma|_{[0,t]}}(\mathbf{J},\mathbf{J})\\
&\leq I_{\gamma|_{[0,t]}}(\mathbf{V}_t,\mathbf{V}_t)\\
&\leq I_c(\overline u_t,\overline u_t)\\
&=\frac{c_c(t)}{s_c(t)}|\mathbf{J}(t)|_{\mathbf{g}}^2.
\end{align*}
Therefore,
\[
\frac{d}{dt}\log\left(\frac{|\mathbf{J}(t)|_{\mathbf{g}}}{s_c(t)}\right)\leq0,
\]
and the same limit at $0$ shows that
\[
|\mathbf{J}(t)|_{\mathbf{g}}\leq s_c(t)|w|_{\mathbf{g}}.
\]
Continuity allows us to include the first conjugate time and the endpoint of the
indicated interval. If $c>0$, $b\geq\rho_c$, and there were no conjugate point
before $\rho_c$, the inequality, applied for $t<\rho_c$, would give
$\mathbf{J}(\rho_c)=0$ for every normal Jacobi field with $\mathbf{J}(0)=0$. Taking
$w\ne0$ would produce a conjugate point at $\rho_c$. Consequently,
$\tau\leq\rho_c$, completing (b).

It remains to study equality. In either case, the quotient
$\frac{|\mathbf{J}|_{\mathbf{g}}}{s_c}$ is monotone and its limit at zero is $|w|_{\mathbf{g}}$. If equality holds
at $t_0$, this quotient is constant on $[0,t_0]$. In the preceding chains of
inequalities, equality must hold both in the comparison of the
index forms and in the curvature estimate. The equality condition
for the model identity in the first case and that of
Lemma~\ref{lem:indice-campos-jacobi} in the second give
\[
\mathbf{J}(t)=s_c(t)P_tw.
\]
Moreover, the self-adjoint endomorphism
$\mathbf{Z}\mapsto\mathbf{R}(\mathbf{Z},\dot\gamma(t))\dot\gamma(t)$
of the subspace $\dot\gamma(t)^\perp\subseteq T_{\gamma(t)}M$
is bounded above by $cI$ in the first case and below by $cI$ in the
second. Since $P_tw\in\dot\gamma(t)^\perp$, equality of the quadratic
forms at $P_tw$ implies
\[
\mathbf{R}(P_tw,\dot\gamma)\dot\gamma=cP_tw.
\]
Conversely, these two identities turn all the preceding steps
into equalities and yield equality in the comparison.
\end{proof}

\subsection{Cut times and injectivity radius}

The absence of conjugate points guarantees local invertibility of
the exponential map but does not ensure global injectivity.
Two different geodesics may reach the same point before its
differential becomes singular. Cut time describes how far a geodesic
realizes distance and distinguishes these two causes of failure
of minimization.

\begin{definition}[Cut time and cut locus]
\label{def:tiempo-corte-lugar-corte}
\index{cut time}\index{cut locus}
Let \((M,\mathbf g)\) be a complete, connected Riemannian manifold
without boundary, of dimension \(n\geq1\). Write
\[
 S_pM:=\{v\in T_pM\mid |v|_{\mathbf g_p}=1\},
 \qquad
 SM:=\{(p,v)\in TM\mid v\in S_pM\}.
\]
A local orthonormal frame identifies \(SM\) with the product of the
open base set and \(\mathbb S^{n-1}\); these trivializations give
its smooth structure and topology.
For \(v\in S_pM\), define
\[
 c(p,v):=
 \displaystyle\sup\{t>0\mid d_{\mathbf g}(p,\exp_p(tv))=t\}\in(0,\infty].
\]
If \(c(p,v)<\infty\), the point \(\exp_p(c(p,v)v)\) is the
\textbf{cut point} of \(p\) in that direction.
The \textbf{cut locus} of \(p\) is
\[
 \operatorname{Cut}(p)
 :=\{\exp_p(c(p,v)v)\mid v\in S_pM,\ c(p,v)<\infty\}.
\]
Completeness ensures that every exponential map in this definition
is defined.
\end{definition}

If the unit-speed geodesic \(\gamma_v(t)=\exp_p(tv)\) minimizes on
\([0,T]\), it also minimizes on each \([0,t]\), with
\(0\leq t\leq T\): a shorter curve to \(\gamma_v(t)\), followed by
the remaining segment, would contradict minimization up to \(T\).
Radial minimization shows that the set of times in the definition
contains an initial interval of positive length.
Continuity of distance also shows that, when \(c(p,v)\) is finite,
the geodesic remains minimizing at the cut time itself.

\begin{proposition}[Alternatives at the cut time]
\label{prop:alternativas-tiempo-corte}
Under the preceding hypotheses, let \(v\in S_pM\) and \(c=c(p,v)\).
\begin{enumerate}[label=(\alph*)]
\item If \(0<t<c\), the geodesic \(\gamma_v|_{[0,t]}\) has no
points conjugate to \(p\) and is the unique unit-speed minimizing
geodesic between its endpoints.
\item If \(c<\infty\), the geodesic minimizes on \([0,c]\), and
at least one of the following holds:
\(\gamma_v(c)\) is conjugate to \(p\) along \(\gamma_v\), or there
is another unit-speed minimizing geodesic from \(p\) to \(\gamma_v(c)\).
\end{enumerate}
\end{proposition}
\begin{proof}
For (a), choose \(T\) with \(t<T<c\).
The geodesic minimizes up to \(T\), so
Theorem~\ref{teo:conjugado-interior-impide-minimizacion} excludes
every conjugate point in \((0,t]\).
If there were another unit-speed minimizing geodesic
\(\sigma\colon[0,t]\longrightarrow M\) with the same endpoints,
its final velocity would differ: if they agreed, uniqueness of
geodesics with prescribed final data would give
\(\sigma=\gamma_v|_{[0,t]}\).
Concatenating \(\sigma\) with \(\gamma_v|_{[t,T]}\) would produce a
curve of length \(T\), hence a minimizing curve, with a corner at
\(t\). This contradicts
Corollary~\ref{cor:curva-minimizante-es-geodesica}.

For (b), minimization on \([0,c]\) follows by continuity as above.
Suppose that \(\gamma_v(c)\) is not conjugate to \(p\).
Take a decreasing sequence of times \(t_j>c\) with \(t_j\to c^{+}\).
By Hopf--Rinow, there is a unit-speed minimizing geodesic from \(p\)
to \(\gamma_v(t_j)\), with initial velocity \(w_j\in S_pM\)
and length
\[
 \ell_j:=d_{\mathbf g}(p,\gamma_v(t_j))<t_j.
\]
Moreover, \(\ell_j\to c\).
Compactness of \(S_pM\) gives a subsequence with
\(w_j\to w\in S_pM\).
From \(\exp_p(\ell_jw_j)=\exp_p(t_jv)\) we obtain
\(\exp_p(cw)=\exp_p(cv)\).
The geodesic with initial velocity \(w\) has length
\(c=d_{\mathbf g}(p,\gamma_v(c))\) and is minimizing.

If \(w=v\), both \(\ell_jw_j\) and \(t_jv\) would converge to
\(cv\). The absence of conjugacy makes \(d(\exp_p)_{cv}\)
invertible. The map \(\exp_p\) is injective in a neighborhood of
\(cv\), so for large \(j\) we would have
\(\ell_jw_j=t_jv\).
This contradicts \(\ell_j<t_j\), since both direction vectors are
unit vectors. Hence \(w\neq v\), giving the second geodesic.
\end{proof}

\begin{theorem}[Continuity of cut time]
\label{teo:continuidad-tiempo-corte}
If \((M,\mathbf g)\) is complete, connected, and without boundary,
of positive dimension, then
\[
 c\colon SM\longrightarrow(0,\infty]
\]
is continuous, with the extended-order topology on the codomain.
\end{theorem}
\begin{proof}
Let \((p_j,v_j)\to(p,v)\) in \(SM\).
Set \(c_j=c(p_j,v_j)\) and \(c_0=c(p,v)\).
Lemma~\ref{lem:vecindades-normales-convexas-locales} and radial
minimization give \(r>0\) such that \(c_j\geq r\) for all
sufficiently large \(j\): near \(p\), the exponential maps are
diffeomorphisms on tangent balls of a common radius.

If \(c_0<\infty\) and \(T>c_0\), then
\(d_{\mathbf g}(p,\exp_p(Tv))<T\).
Continuity of the exponential map and of distance implies
\(d_{\mathbf g}(p_j,\exp_{p_j}(Tv_j))<T\) for large \(j\).
Thus \(c_j<T\), and letting \(T\to c_0^{+}\) gives
\[
 \limsup_{j\to\infty}c_j\leq c_0.
\]
If \(c_0=\infty\), this inequality imposes no additional condition.

For the reverse inequality, suppose that
\(b:=\liminf_{j\to\infty}c_j<c_0\).
Pass to a subsequence with \(c_j\to b<\infty\), where
\(b\geq r>0\), and every \(c_j\) is finite. Let
\[
 q_j:=\exp_{p_j}(c_jv_j),\qquad q:=\exp_p(bv).
\]
Then \(q_j\to q\) and
\(d_{\mathbf g}(p,q)=\lim_j d_{\mathbf g}(p_j,q_j)=b\).
Since \(b<c_0\), the preceding proposition shows that \(q\) is
not conjugate to \(p\) along \(\gamma_v\), and that
\(\gamma_v|_{[0,b]}\) is the unique unit-speed minimizer.

If there are infinitely many \(j\) for which \(q_j\) is conjugate
to \(p_j\) in direction \(v_j\), the vertical differentials
\(d(\exp_{p_j})_{c_jv_j}\) are singular.
In local trivializations, their matrices converge to that of
\(d(\exp_p)_{bv}\).
Continuity of the determinant would make the latter singular,
a contradiction.

In the remaining case, for every sufficiently large \(j\) there is
\(w_j\in S_{p_j}M\), distinct from \(v_j\), such that
\[
 \exp_{p_j}(c_jw_j)=q_j.
\]
In a local orthonormal frame
\((\mathbf E_1,\ldots,\mathbf E_n)\) near \(p\), write
\[
 w_j=\displaystyle\sum_{i=1}^n a_j^i\mathbf E_i(p_j),
 \qquad (a_j^1,\ldots,a_j^n)\in\mathbb S^{n-1}.
\]
A subsequence of these coefficients converges to
\((a^1,\ldots,a^n)\in\mathbb S^{n-1}\).
For \(w=\displaystyle\sum_{i=1}^n a^i\mathbf E_i(p)\), we have
\(\exp_p(bw)=q\).
This geodesic also has length \(b\), so uniqueness gives \(w=v\).

Now consider \(\mathcal A(x,u)=(x,\exp_x(u))\) near \((p,bv)\).
Its differential is block triangular, with the identity in the
base block and \(d(\exp_p)_{bv}\) in the vertical block.
It is invertible, so \(\mathcal A\) is injective on a neighborhood
of \((p,bv)\).
For large \(j\), both \((p_j,c_jv_j)\) and \((p_j,c_jw_j)\) lie
in that neighborhood and have the same image \((p_j,q_j)\).
This contradicts \(w_j\neq v_j\).

We have proved \(c_0\leq\liminf_jc_j\).
Together with the upper bound, this gives \(c_j\to c_0\).
When \(c_0=\infty\), the same argument excludes every bounded
subsequence of \(c_j\), which is precisely convergence to infinity.
\end{proof}

\begin{proposition}[Injectivity radius and cut times]
\label{prop:inyectividad-infimo-tiempos-corte}
If \((M,\mathbf g)\) is complete, connected, and without boundary,
of positive dimension, then
\[
 \operatorname{inj}(p)
 =\displaystyle\inf_{v\in S_pM}c(p,v)
 =d_{\mathbf g}(p,\operatorname{Cut}(p)).
\]
Distance to the empty set is interpreted as \(+\infty\).
\end{proposition}
\begin{proof}
Set \(a=\displaystyle\inf_{v\in S_pM}c(p,v)\).
If \(0<r<a\), every vector \(u\) with
\(0<|u|_{\mathbf g_p}<r\) lies strictly before the cut time in its
direction.
Proposition~\ref{prop:alternativas-tiempo-corte} and the
characterization of conjugate points show that \(d(\exp_p)_u\)
is invertible. It is also invertible at \(u=0\).
Moreover, if \(\exp_p(u)=\exp_p(w)\) with \(|u|_{\mathbf g_p},|w|_{\mathbf g_p}<r\),
both radial geodesics minimize and have the same length.
If that length is zero, \(u=w=0\); otherwise, uniqueness before
the cut time gives \(u=w\).
Thus \(\exp_p\) is an injective local diffeomorphism on
\(B_{\mathbf g_p}(0,r)\), hence a diffeomorphism onto its image.
Definition~\ref{def: radio de inyectividad} yields
\(\operatorname{inj}(p)\geq a\), also when \(a=\infty\).

If \(r>a\), there is \(v\in S_pM\) with \(c(p,v)<r\).
At the vector \(c(p,v)v\), either the differential of \(\exp_p\)
is singular, or there is another vector of the same norm with
the same image, by the alternative at the cut time.
In either case, \(\exp_p\) is not a diffeomorphism on the ball
of radius \(r\). Hence \(\operatorname{inj}(p)\leq a\).
Finally, the cut point in direction \(v\) has distance \(c(p,v)\)
from \(p\), proving the second equality.
\end{proof}

\begin{proposition}[The injectivity domain]
\label{prop:dominio-inyectividad-lugar-corte}
Under the same hypotheses, let
\[
 \mathcal I_p:=\{0\}\cup
 \left\{u\in T_pM\setminus\{0\}\ \middle|\
 |u|_{\mathbf g_p}<c\!\left(p,\frac{u}{|u|_{\mathbf g_p}}\right)\right\}.
\]
Then \(\mathcal I_p\) is open and star-shaped about zero,
\(\operatorname{Cut}(p)\) is closed, and
\[
 \exp_p\colon\mathcal I_p\longrightarrow M\setminus\operatorname{Cut}(p)
\]
is a diffeomorphism.
\end{proposition}
\begin{proof}
Continuity of cut time shows that \(\mathcal I_p\) is open away
from zero. Radial minimization and the alternative at the cut time
show that it contains a ball about zero, so it is open there as well.
It is star-shaped by definition.
The argument in the preceding proposition shows that \(\exp_p\)
is an injective local diffeomorphism on \(\mathcal I_p\).

By Hopf--Rinow, every \(q\neq p\) is the endpoint of a unit-speed
minimizing geodesic from \(p\), of length
\(t=d_{\mathbf g}(p,q)\leq c(p,v)\).
If \(q\) is not in the cut locus, this inequality is strict and
\(tv\in\mathcal I_p\).
Conversely, a point reached strictly before the cut time cannot
be a cut point in another direction: both geodesics would minimize
and have the same length, so uniqueness before the cut time would
make them identical.
Nor can \(p\) be a cut point, since every cut time is positive.
Thus the image is exactly \(M\setminus\operatorname{Cut}(p)\).

To prove that the cut locus is closed, let \(q_j\to q\), where
\(q_j=\exp_p(c(p,v_j)v_j)\) and \(v_j\in S_pM\).
Passing to a subsequence gives \(v_j\to v\in S_pM\), and
\[
 c(p,v_j)=d_{\mathbf g}(p,q_j)\longrightarrow d_{\mathbf g}(p,q)<\infty.
\]
Continuity of cut time gives \(c(p,v)=d_{\mathbf g}(p,q)\);
continuity of the exponential map then yields
\(q=\exp_p(c(p,v)v)\).
Hence \(q\in\operatorname{Cut}(p)\).
\end{proof}

\begin{proposition}[Continuity of the injectivity radius]
\label{prop: continuidad radio inyectividad conexa y completa}
\index{injectivity radius!continuity}
Let \((M,\mathbf g)\) be a complete, connected Riemannian manifold
without boundary. Then
\[
 \operatorname{inj}\colon M\longrightarrow(0,\infty]
\]
is continuous, where \((0,\infty]\) has the extended-order topology.
\end{proposition}
\begin{proof}
If \(\dim M=0\), connectedness makes \(M\) a single point and
\(\operatorname{inj}\equiv\infty\). Assume \(\dim M=n\geq1\).
By Proposition~\ref{prop:inyectividad-infimo-tiempos-corte},
\[
 \operatorname{inj}(x)=\displaystyle\inf_{v\in S_xM}c(x,v).
\]
Fix \(p\in M\), a sequence \(p_j\to p\), and a local orthonormal
frame \((\mathbf E_1,\ldots,\mathbf E_n)\) near \(p\).
For \(x\) in this open set and \(\theta\in\mathbb S^{n-1}\), put
\[
 v(x,\theta):=\displaystyle\sum_{i=1}^n\theta^i\mathbf E_i(x).
\]
The function \((x,\theta)\mapsto c(x,v(x,\theta))\) is continuous.
Its infimum on each fiber is attained: if it is finite, take a
minimizing sequence and a convergent subsequence in the compact
sphere; if it is infinite, every value is infinite.

If \(\operatorname{inj}(p)<\infty\), choose \(\theta_0\) where
this value is attained. Then
\[
 \operatorname{inj}(p_j)\leq c(p_j,v(p_j,\theta_0))
 \longrightarrow c(p,v(p,\theta_0))=\operatorname{inj}(p),
\]
proving the upper-limit inequality.
For the lower limit, if
\(b:=\liminf_j\operatorname{inj}(p_j)<\infty\), take a
subsequence along which these radii converge to \(b\), together
with minimizers \(\theta_j\in\mathbb S^{n-1}\).
A further subsequence satisfies \(\theta_j\to\theta\), and
continuity gives
\[
 b=\lim_j c(p_j,v(p_j,\theta_j))
 =c(p,v(p,\theta))\geq\operatorname{inj}(p).
\]
If the lower limit is infinite, the inequality is immediate.
When \(\operatorname{inj}(p)<\infty\), the two bounds prove
convergence to its value at \(p\).
When \(\operatorname{inj}(p)=\infty\), the second bound excludes
a finite lower limit and proves convergence to infinity.
\end{proof}
\section{Integration on Riemannian manifolds}

Integration on a manifold requires replacing the Euclidean volume element by an object that transforms correctly under changes of coordinates. On an oriented manifold, differential forms of top degree provide such an object; the Riemannian metric will later select a canonical volume form. We begin with orientation and integration of forms in charts.
\begin{definition}\label{def:variedades-riemannianas-consistentemente-orientadas}\index{consistently oriented}
 Let $V$ be a finite-dimensional vector space. We say that two ordered bases $(E_{1},\dots,E_{n})$ and $(\widetilde{E}_{1},\dots,\widetilde{E}_{n})$ are \textbf{consistently oriented} if the change-of-basis matrix $(B_{i}^{j})$ satisfying $E_{i}=B_{i}^{j}\widetilde{E}_{j}$ has positive determinant.
\end{definition}
One can show that consistent orientation is an equivalence relation on the set of ordered bases of a finite-dimensional vector space. In positive dimension this relation has exactly two classes: the sign of the determinant relative to a fixed basis distinguishes them, and reversing one basis vector interchanges them. This gives the following definition:
\begin{definition}\label{def:variedades-riemannianas-orientacion-para}\index{orientation}
 Let $V$ be a positive-dimensional finite-dimensional vector space. An \textbf{orientation of $V$} is an equivalence class of ordered bases. Given an ordered basis $(E_{1},\dots,E_{n})$, denote the orientation it determines by $[E_{1},\dots,E_{n}]$. An \textbf{oriented vector space} is a vector space $V$ together with an orientation determined by an ordered basis $[E_{1},\dots,E_{n}]$. Once we fix an orientation of $V$, a basis with this orientation is said to be \textbf{positively oriented}. An ordered basis that is not positively oriented is said to be \textbf{negatively oriented}.
\end{definition}

In dimension zero, we assign an orientation sign $+1$ or $-1$ to each point by convention. This permits orientation and integration on the boundary of a one-dimensional manifold.

\begin{definition}\label{def:variedades-riemannianas-orientacion-puntual}\index{pointwise orientation}
 Let $M$ be a smooth manifold with or without boundary. A \textbf{pointwise orientation} on $M$ is a choice of orientation on each tangent space. If $(\mathbf{E}_{i})$ is a local frame of $TM$, we say that $(\mathbf{E}_{i})$ \textbf{is positively oriented} if $(\mathbf{E}_{1}|_{p},\dots,\mathbf{E}_{n}|_{p})$ is a positively oriented ordered basis of $T_{p}M$ for each $p\in U$, where $U$ is the domain of the local frame $(\mathbf{E}_{i})$. If $(\mathbf{E}_{1}|_{p},\dots,\mathbf{E}_{n}|_{p})$ is negatively oriented for each $p\in U$, we say that the local frame is \textbf{negatively oriented}.

 A \textbf{continuous orientation of $M$} is a pointwise orientation such that, for each $p\in M$, there exists an open set $U\subseteq M$ such that $U$ is the domain of a positively oriented local frame. We say that $M$ is \textbf{orientable} if it has a continuous pointwise orientation. We denote this orientation by $\mathscr{O}$ and say that $(M,\mathscr{O})$ is an \textbf{oriented smooth manifold}.

 Finally, if $(U,\phi)$ is a smooth chart of an oriented manifold $M$, with $\phi=(x^{1},\dots,x^{n})$, we say that it is a \textbf{positively oriented chart} (or simply an \textbf{oriented chart}) if the coordinate frame $\bigl(\boldsymbol{\partial}_i\bigr)$ is positively oriented. If this coordinate frame is negatively oriented, we say that the chart is a \textbf{negatively oriented chart}.
\end{definition}
\begin{proposition}[The Riemannian volume form]\label{forma de volumen}\index{volume form!Riemannian}
 Let $(M,\mathbf{g})$ be an oriented Riemannian manifold of dimension $n\geq1$ with or without boundary. There exists a unique $n$-form $dV_{\mathbf{g}}\in \Omega^{n}(M)$, called the \textbf{Riemannian volume form}, characterized by any of the following properties:
 \begin{enumerate}[label=(\alph*)]
 \item If $(\boldsymbol{\varepsilon}^{1},\dots,\boldsymbol{\varepsilon}^{n})$ is any positively oriented orthonormal local coframe of $T^{*}M$, then \[dV_{\mathbf{g}}=\boldsymbol{\varepsilon}^{1}\wedge\cdots\wedge\boldsymbol{\varepsilon}^{n}.\]
 \item If $(\mathbf{E}_{1},\dots,\mathbf{E}_{n})$ is any positively oriented orthonormal local frame of $TM$, then \[dV_{\mathbf{g}}(\mathbf{E}_{1},\dots,\mathbf{E}_{n})=1.\]
 \item If $(U,\phi)$, $\phi=(x^{1},\dots,x^{n})$ is an oriented chart of $M$, then \[dV_{\mathbf{g}}=\sqrt{\det(\mathbf{g})}\,\mathbf{d}x^{1}\wedge\cdots\wedge \mathbf{d}x^{n}.\]
 \end{enumerate}
 \end{proposition}
\begin{proof}
Let $(U,x^1,\dots,x^n)$ be a positively oriented chart. The matrix
$G_x=(g_{ij})$ is symmetric and positive definite, so
$\det G_x>0$. On $U$, define
\begin{equation}
\label{eq:construccion-local-forma-volumen}
\boldsymbol{\omega}_x:=\sqrt{\det G_x}\,
\mathbf{d}x^1\wedge\cdots\wedge \mathbf{d}x^n.
\end{equation}
This form is smooth. We will show that the forms obtained in this way agree on
overlaps.

Let $(x^1,\dots,x^n)$ and $(y^1,\dots,y^n)$ be two positively oriented coordinate
systems and write
$G_y=(\widetilde g_{ab})$ for the matrix of the metric in the coordinates
$y$. On the overlap,
\[
\widetilde g_{ab}
=g_{ij}\frac{\partial x^i}{\partial y^a}
          \frac{\partial x^j}{\partial y^b},
\]
that is,
\[
G_y=
\left(\frac{\partial x}{\partial y}\right)^{\!T}
G_x
\left(\frac{\partial x}{\partial y}\right).
\]
Taking determinants gives
\begin{equation}
\label{eq:transformacion-determinante-metrica-volumen}
\det G_y
=\det\left(\frac{\partial x}{\partial y}\right)^2
\det G_x.
\end{equation}
Both coordinate systems induce the fixed orientation; therefore,
\[
\det\left(\frac{\partial x}{\partial y}\right)>0.
\]
The transformation law for the wedge product, in turn, gives
\[
\mathbf{d}y^1\wedge\cdots\wedge \mathbf{d}y^n
=\det\left(\frac{\partial y}{\partial x}\right)
\mathbf{d}x^1\wedge\cdots\wedge \mathbf{d}x^n.
\]
Combining this equality with
\eqref{eq:transformacion-determinante-metrica-volumen}, we conclude that
\[
\begin{aligned}
\sqrt{\det G_y}\,\mathbf{d}y^1\wedge\cdots\wedge \mathbf{d}y^n
&=\det\left(\frac{\partial x}{\partial y}\right)
  \sqrt{\det G_x}
  \det\left(\frac{\partial y}{\partial x}\right)
  \mathbf{d}x^1\wedge\cdots\wedge \mathbf{d}x^n\\
&=\sqrt{\det G_x}\,\mathbf{d}x^1\wedge\cdots\wedge \mathbf{d}x^n.
\end{aligned}
\]
Thus the local forms glue together to define a smooth global form
$dV_{\mathbf{g}}$. Property (c) holds by construction.

Now let $(\mathbf{E}_1,\dots,\mathbf{E}_n)$ be a positively oriented orthonormal local
frame. In an oriented chart, write
\[
\mathbf{E}_a=B_a^{i}\boldsymbol{\partial}_i
\]
and denote by $B$ the matrix whose columns are the components of the
$\mathbf{E}_a$. Orthonormality is equivalent to
\[
B^T G_xB=I.
\]
Therefore, $(\det B)^2\det G_x=1$. Since the frame and the chart have the
same orientation, $\det B>0$, and hence
\[
\det B=\frac{1}{\sqrt{\det G_x}}.
\]
It follows from \eqref{eq:construccion-local-forma-volumen} that
\[
dV_{\mathbf{g}}(\mathbf{E}_1,\dots,\mathbf{E}_n)
=\sqrt{\det G_x}\,\det B=1.
\]
This proves (b). If $(\boldsymbol{\varepsilon}^1,\dots,\boldsymbol{\varepsilon}^n)$ is the dual
coframe, then
$(\boldsymbol{\varepsilon}^1\wedge\cdots\wedge\boldsymbol{\varepsilon}^n)(\mathbf{E}_1,\dots,\mathbf{E}_n)=1$.
Since the space of $n$-covectors at each point is one-dimensional, (b)
implies
\[
dV_{\mathbf{g}}=\boldsymbol{\varepsilon}^1\wedge\cdots\wedge\boldsymbol{\varepsilon}^n,
\]
which is (a). Conversely, (a) immediately implies (b), and the calculation
with the coordinate frame shows that (b) implies (c). Thus the three
characterizations are equivalent.

The construction does not depend on the orthonormal frame. Two
positively oriented orthonormal frames are related by a matrix
with values in $SO(n)$, whose determinant is one; the wedge products of their dual
coframes agree. If another $n$-form satisfies any of the
characterizations, its value on each positively oriented orthonormal
frame is one, so it agrees pointwise with $dV_{\mathbf{g}}$. This
proves uniqueness.

\end{proof}

\begin{remark}\label{obs:distincion-forma-densidad-medida-riemanniana}
The notation $dV_{\mathbf{g}}$ here denotes an $n$-form and requires a fixed
orientation. The Riemannian density $\boldsymbol{\mu}_{\mathbf{g}}$, defined
later, exists even when $M$ is not orientable and transforms with the absolute
value of the Jacobian. The Riemann--Lebesgue measure $\lambda_{\mathbf{g}}$ will be
constructed by integrating this density; the notation $d\lambda_{\mathbf{g}}$ indicates only
integration with respect to the measure $\lambda_{\mathbf{g}}$ and does not represent a differential
form or an exterior derivative.
\end{remark}
The volume form has the important property of compatibility with the Levi--Civita connection. To prove this, we first recall an important result from differential calculus:
\begin{proposition}[Jacobi's formula]\label{formula de jacobi}\index{Jacobi's formula}
If $A\colon I\subseteq \mathbb{R}\longrightarrow GL(n,\mathbb{R})$ is differentiable, then
\[
\displaystyle\frac{d}{dt}\det\big(A(t)\big)
= \det\big(A(t)\big)\operatorname{tr}\left(A(t)^{-1}\displaystyle\frac{d}{dt}A(t)\right).
\]

\end{proposition}
\begin{proof}
Fix \(t\in I\) and denote the columns of \(A(t)\) by
\(\mathbf a_1(t),\dots,\mathbf a_n(t)\). Multilinearity
of the determinant gives
\[
\frac d{dt}\det A(t)
=\displaystyle\sum_{j=1}^n
\det\bigl(\mathbf a_1(t),\dots,\mathbf a'_j(t),\dots,\mathbf a_n(t)\bigr).
\]
Set \(B(t):=A(t)^{-1}A'(t)\). Then \(A'=AB\), and
\[
\mathbf a'_j=\displaystyle\sum_{i=1}^nB_{ij}\mathbf a_i.
\]
Substituting this sum into the determinant indexed by \(j\), all the
terms with \(i\ne j\) vanish because they contain two equal columns.
The remaining term is \(B_{jj}\det A\). Therefore,
\[
\frac d{dt}\det A(t)
=\det A(t)\displaystyle\sum_{j=1}^nB_{jj}(t)
=\det A(t)\operatorname{tr}\bigl(A(t)^{-1}A'(t)\bigr).
\]
\end{proof}
\begin{remark}\label{derivada raiz del determinante de la metrica}
Proposition~\ref{formula de jacobi} allows us to compute the partial
derivatives of $\sqrt{\det(\mathbf{g})}$ and relate them to the Christoffel
symbols. Consider a chart $(U,\phi)$ with induced coordinates
\[
\phi=(x^{1},\dots,x^{n})\colon U\longrightarrow \mathbb{R}^{n}.
\]
Denote by $G=(g_{ij})$ the matrix of metric components, where
\[
g_{ij}(p)=\mathbf{g}_{p}\left(\boldsymbol{\partial}_i\big|_{p},\boldsymbol{\partial}_j\big|_{p}\right),
\qquad p\in U.
\]

Fix a point $p\in U$ and a coordinate $x^{a}$. Define a curve in $U$ by
\[
\gamma(t)=\phi^{-1}(x^{1}(p),\dots,x^{a}(p)+t,\dots,x^{n}(p)).
\]
In coordinates, we have
\[
(\phi\circ \gamma)(t)=(x^{1}(p),\dots,x^{a}(p)+t,\dots,x^{n}(p)).
\]

Let $f\colon U\longrightarrow \mathbb{R}$ be a smooth function and write $F=f\circ \phi^{-1}\colon \phi(U)\subseteq \mathbb{R}^{n}\longrightarrow \mathbb{R}$. Then
\[
f(\gamma(t))=F((\phi\circ \gamma)(t)).
\]
Applying the chain rule,
\[
\displaystyle\frac{d}{dt}\bigg|_{t=0}f(\gamma(t))
=dF_{\phi(p)}\left(\displaystyle\frac{d}{dt}\bigg|_{t=0}(\phi\circ \gamma)(t)\right).
\]
Since
\[
\displaystyle\frac{d}{dt}\bigg|_{t=0}(\phi\circ \gamma)(t)=(0,\dots,1,\dots,0)=e_{a}\in \mathbb{R}^{n},
\]
it follows that
\[
\displaystyle\frac{d}{dt}\bigg|_{t=0}f(\gamma(t))
=dF_{\phi(p)}(e_{a})
=\displaystyle\frac{\partial F}{\partial u^{a}}(\phi(p)),
\]
where $(u^{1},\dots,u^{n})$ are the standard coordinates of $\mathbb{R}^{n}$. By definition, this agrees with the action of $\displaystyle\boldsymbol{\partial}_{a}|_{p}$ on $f$. Consequently,
\[
\gamma'(0)=\boldsymbol{\partial}_a\big|_{p}.
\]

Applying this to each coordinate component of the metric, which is a smooth function $g_{ij}\colon U\longrightarrow \mathbb{R}$, we conclude that
\[
\displaystyle\frac{d}{dt}\bigg|_{t=0}g_{ij}(\gamma(t))=\partial_{a}g_{ij}(p).
\]

Defining the curve of matrices
\[
A(t)=(g_{ij}(\gamma(t))),
\]
and applying Proposition~\ref{formula de jacobi} gives
\[
\displaystyle\frac{d}{dt}\det\big(A(t)\big)
= \det\big(A(t)\big)\operatorname{tr}\left(A(t)^{-1}\displaystyle\frac{d}{dt}A(t)\right).
\]

Evaluating at $t=0$ and using $\gamma(0)=p$ together with the preceding calculation, we conclude that
\[
\partial_{a}(\det G)(p)=(\det G(p))g^{ij}(p)\partial_{a}g_{ij}(p).
\]

Thus Jacobi's formula in one real variable becomes, in local coordinates on the manifold, the identity
\[
\partial_{a}(\det G)=(\det G)g^{ij}\partial_{a}g_{ij}.
\] It follows that
\[
\partial_{a}\sqrt{\det G}
=\displaystyle\frac{1}{2\sqrt{\det G}}\partial_{a}(\det G).
\]
Substituting the preceding identity gives
\[
\partial_{a}\sqrt{\det G}
=\displaystyle\frac{1}{2\sqrt{\det G}}\big((\det G)g^{ij}\partial_{a}g_{ij}\big),
\]
and therefore
\[
\partial_{a}\sqrt{\det G}
=\displaystyle\frac{1}{2}\sqrt{\det G}\,g^{ij}\partial_{a}g_{ij}.
\]

In particular, using Corollary~\ref{expresion local christoffel levi civita}, this formula gives the identity
\[
\Gamma^{i}_{ik}
=\displaystyle\frac{1}{\sqrt{\det G}}\partial_{k}\left(\sqrt{\det G}\right),
\]
where $\Gamma^{i}_{jk}$ denote the Christoffel symbols of the Levi--Civita connection.
\end{remark}

\begin{proposition}\label{forma de volumen paralela}
Let $(M,\mathbf{g})$ be an oriented Riemannian manifold with or without boundary and with Levi--Civita connection $\nabla$. Then the Riemannian volume form is parallel with respect to the connection, that is,
\[
\nabla dV_{\mathbf{g}}=0.
\]
\end{proposition}

\begin{proof}
In an oriented chart $(x^{1},\dots,x^{n})$, we have
\[
dV_{\mathbf{g}}=\sqrt{\det G}\,\mathbf{d}x^{1}\wedge\cdots\wedge \mathbf{d}x^{n}.
\]
Using Proposition~\ref{conexion producto cuña} and linearity of the connection gives
\[
\nabla_{a}dV_{\mathbf{g}}
=\nabla_{a}(\sqrt{\det G})\,\mathbf{d}x^{1}\wedge\cdots\wedge \mathbf{d}x^{n}
+\sqrt{\det G}\,\nabla_{a}(\mathbf{d}x^{1}\wedge\cdots\wedge \mathbf{d}x^{n}).
\]

For the first term, since $\sqrt{\det G}$ is a function, we have
\[
\nabla_{a}(\sqrt{\det G})=\partial_{a}(\sqrt{\det G}).
\]

For the second term, applying Proposition~\ref{conexion producto cuña} again,
\[
\nabla_{a}(\mathbf{d}x^{1}\wedge\cdots\wedge \mathbf{d}x^{n})
=\displaystyle\sum_{s=1}^{n}\mathbf{d}x^{1}\wedge\cdots\wedge \nabla_{a}\mathbf{d}x^{s}\wedge\cdots\wedge \mathbf{d}x^{n}.
\]
In coordinates,
\[
\nabla_{a}\mathbf{d}x^{s}=-\Gamma^{s}_{ar}\mathbf{d}x^{r}.
\]
If $r\neq s$, the term
\[
\mathbf{d}x^{1}\wedge\cdots\wedge(-\Gamma^{s}_{ar}\mathbf{d}x^{r})\wedge\cdots\wedge \mathbf{d}x^{n}
\]
contains two equal factors ($\mathbf{d}x^{r}$ appears twice), so it vanishes. The only nontrivial contribution occurs when $r=s$, giving
\[
\nabla_{a}(\mathbf{d}x^{1}\wedge\cdots\wedge \mathbf{d}x^{n})
=-\Big(\Gamma^{1}_{a1}+\cdots+\Gamma^{n}_{an}\Big)\mathbf{d}x^{1}\wedge\cdots\wedge \mathbf{d}x^{n}.
\]

Combining both results,
\[
\nabla_{a}dV_{\mathbf{g}}
=\Big(\partial_{a}(\sqrt{\det G})-\sqrt{\det G}(\Gamma^{1}_{a1}+\cdots+\Gamma^{n}_{an})\Big)\mathbf{d}x^{1}\wedge\cdots\wedge \mathbf{d}x^{n}.
\]

By Remark~\ref{derivada raiz del determinante de la metrica},
\[
\partial_{a}(\sqrt{\det G})=\displaystyle\frac{1}{2}\sqrt{\det G}\,g^{ij}\partial_{a}g_{ij}.
\]
Since the Levi--Civita connection is compatible with $\mathbf{g}$, we have
\[
\Gamma^{r}_{ra}=\displaystyle\frac{1}{2}g^{ij}\partial_{a}g_{ij}.
\]
Moreover, symmetry in the lower indices implies $\Gamma^{r}_{ar}=\Gamma^{r}_{ra}$, so
\[
\Gamma^{1}_{a1}+\cdots+\Gamma^{n}_{an}=\Gamma^{r}_{ra}.
\]

Substituting these identities,
\[
\nabla_{a}dV_{\mathbf{g}}
=\big(\Gamma^{r}_{ra}\sqrt{\det G}-\Gamma^{r}_{ra}\sqrt{\det G}\big)\mathbf{d}x^{1}\wedge\cdots\wedge \mathbf{d}x^{n}=0.
\]

Thus, $dV_{\mathbf{g}}$ is parallel with respect to the Levi--Civita connection.
\end{proof}

 We now have the ingredients needed to define integrals of $n$-forms on smooth manifolds. To do so, we use Theorem~\ref{particionesdelaunidad}.
 \begin{definition}
\label{def:variedades-riemannianas-integral-de-sobre}
\index{integration of forms on Euclidean domains}
Let $D\subseteq\mathbb R^n$, or $D\subseteq\mathbb H^n$, be
Jordan measurable and let $\boldsymbol\omega$ be a continuous
$n$-form on $\overline D$. If
$\boldsymbol\omega=f\,\mathbf d x^1\wedge\cdots\wedge\mathbf d x^n$,
define
\[
\int_D\boldsymbol\omega:=\int_D f\,dx^1\cdots dx^n.
\]
For an open subset $U$ of $\mathbb R^n$, or a relatively open
subset of $\mathbb H^n$, and a continuous $n$-form with compact
support $K\subseteq U$, define
\[
\int_U\boldsymbol\omega:=\int_D\boldsymbol\omega,
\qquad K\subseteq D,\quad\overline D\subseteq U,
\]
where $D$ is Jordan measurable. Such a $D$ exists: around each
point of $K$, choose a sufficiently small closed box contained in
$U$, intersected with $\mathbb H^n$ when appropriate, whose relative
interior contains that point. A finite subcover of these interiors
covers $K$, and the union of the corresponding boxes serves as $D$.
Independence of this choice is proved below.
\end{definition}
 \begin{definition}\label{def:variedades-riemannianas-integral-variedad-suave-orientada-frontera}\index{local integral of forms on oriented manifolds}
 Let $M$ be an oriented smooth manifold with or without boundary and let $\boldsymbol{\omega}$ be a continuous $n$-form on $M$ with compact support contained in the domain of a smooth chart $(U,\phi)$ that is positively or negatively oriented. Define the integral of $\boldsymbol{\omega}$ over $M$ to be \[\int_{M}\boldsymbol{\omega}=\pm \int_{\phi(U)}(\phi^{-1})^{*}\boldsymbol{\omega},\] where the sign is positive if the chart is positively oriented and negative if it is negatively oriented.
 \end{definition}
 \begin{definition}
\label{def:variedades-riemannianas-integral-variedad-suave-orientada-frontera-2}
\index{integration of forms on oriented manifolds}
Let $M$ be an oriented smooth manifold with or without boundary,
and let $\boldsymbol\omega$ be a continuous $n$-form with compact
support $K$. Choose oriented chart domains $U_1,\ldots,U_k$ covering
$K$, and nonnegative smooth functions $\psi_1,\ldots,\psi_k$ such
that $\supp(\psi_i)\subseteq U_i$ and
$\displaystyle\sum_{i=1}^k\psi_i=1$ near $K$. Obtain them from
a partition of unity subordinate to the cover
$U_1,\ldots,U_k,M\setminus K$, omitting the function associated
with $M\setminus K$, whose support is disjoint from $K$. Define
\[
\int_M\boldsymbol\omega
:=\displaystyle\sum_{i=1}^k\int_M\psi_i\boldsymbol\omega,
\]
where each summand is computed by the preceding local definition.
\end{definition}

\begin{proposition}[Well-definedness of integration of forms]
\label{prop:integral-formas-bien-definida}
The integral is independent of the Jordan measurable set, charts,
finite cover, and partition functions chosen.
\end{proposition}
\begin{proof}
On a Euclidean domain, extend the coefficient $f$ by zero outside
$U$, within $\mathbb R^n$ or $\mathbb H^n$ as appropriate.
The extension is continuous because its support is compact in $U$.
If $D_1,D_2$ satisfy the conditions in the definition, $f$ vanishes
on the difference between their union and either set. Additivity
of the Riemann integral gives
\[
\int_{D_1}f\,dx=\int_{D_1\cup D_2}f\,dx
=\int_{D_2}f\,dx.
\]

For independence of the chart, suppose a form has support contained
in two domains $U,V$, with charts $\phi,\psi$. Let
$\epsilon_\phi,\epsilon_\psi\in\{1,-1\}$ be their orientation
signs. On the overlap, write
\[
\Theta:=\psi\circ\phi^{-1},\qquad
(\phi^{-1})^*\boldsymbol\omega=f_\phi\,\mathbf d x^1\wedge\cdots\wedge\mathbf d x^n,
\qquad
(\psi^{-1})^*\boldsymbol\omega=f_\psi\,\mathbf d y^1\wedge\cdots\wedge\mathbf d y^n.
\]
The transformation law for forms and the orientation give
\[
f_\phi(x)=f_\psi(\Theta(x))\det D\Theta(x),
\qquad
\operatorname{sgn}(\det D\Theta)=\epsilon_\phi\epsilon_\psi.
\]
By Theorem~\ref{teorema de cambio de variable},
\[
\epsilon_\phi\int_{\phi(U)}f_\phi\,dx
=\epsilon_\psi\int_{\psi(V)}f_\psi\,dy.
\]
The integrands vanish outside the images of the support, so the
overlap suffices for the change of variables. In boundary charts,
apply it on the interiors; the boundary hyperplanes have zero
$n$-dimensional measure and do not change the integrals.

Finally, let $(\psi_i)_{i=1}^k$ and $(\eta_j)_{j=1}^{\ell}$ be two
choices. Each $\psi_i\eta_j\boldsymbol\omega$ has compact support
in the intersection of the two relevant charts, so its integral can
be computed in either. Local linearity and the sums equal to one
on the support yield
\[
\begin{aligned}
\displaystyle\sum_{i=1}^k\int_M\psi_i\boldsymbol\omega
&=\displaystyle\sum_{i=1}^k\displaystyle\sum_{j=1}^{\ell}
\int_M\psi_i\eta_j\boldsymbol\omega\\
&=\displaystyle\sum_{j=1}^{\ell}\int_M\eta_j\boldsymbol\omega.
\end{aligned}
\]
This proves independence of both the cover and the partition.
\end{proof}

\medskip

We extend this definition to Riemannian manifolds as follows:

 \begin{definition}\label{def:variedades-riemannianas-integral-variedad-riemanniana-orientada-continua}\index{integral with respect to the volume form}
 Let $(M,\mathbf{g})$ be an oriented Riemannian manifold with or without boundary and let $f\colon M\longrightarrow \mathbb{R}$ be continuous with compact support. Define the integral of $f$ over $M$ to be \[\int_{M}f=\int_{M}fdV_{\mathbf{g}}.\]
 \end{definition}
 \medskip

 The nonorientable case is handled by densities, which incorporate the absolute value of the Jacobian and require no choice of orientation.
 \begin{definition}[Density of weight $w$ on a vector space]\label{def:variedades-riemannianas-densidad-de-peso-w-en-un-espacio-vectorial}\index{density!of weight \(w\) on a vector space}
Let $V$ be a vector space of dimension $n$ and let $w \in \mathbb{R}$.
Write
\[
\operatorname{Fr}(V)
:=
\{\,(v_1,\dots,v_n)\in V^n\mid (v_1,\dots,v_n)
   \text{ is an ordered basis of }V\,\}.
\]
A \textbf{density of weight $w$ on $V$} is a function
\[
\mu\colon \operatorname{Fr}(V)\longrightarrow\mathbb{R}
\]
such that, for every invertible linear operator
$T\colon V\longrightarrow V$ and every
$(v_1,\dots,v_n)\in\operatorname{Fr}(V)$,
\[
\mu(T(v_{1}),\dots,T(v_{n}))=|\det(T)|^{w}\mu(v_{1},\dots,v_{n}).
\]
A density $\mu$ is called \textbf{positive} if
$\mu(v_{1},\dots,v_{n})>0$ for every ordered basis
$(v_{1},\dots,v_{n})$ of $V$.

Denote the set of densities of weight $w$ on $V$ by $\mathcal{D}_{w}(V)$.
\end{definition}

With pointwise operations, \(\mathcal D_w(V)\) is a one-dimensional vector
space. Indeed, fixing an ordered basis
\(\mathcal E=(e_1,\dots,e_n)\), every ordered basis
\(\mathcal V=(v_1,\dots,v_n)\) can be written uniquely as
\(\mathcal V=T\mathcal E\) with \(T\in\operatorname{GL}(V)\). Therefore,
\[
\mu(\mathcal V)=|\det T|^w\mu(\mathcal E).
\]
Thus the value of \(\mu\) at \(\mathcal E\) determines all its values.
Conversely, for each \(c\in\mathbb R\), the formula
\(\mu_c(T\mathcal E):=|\det T|^wc\) is well-defined and determines a
density. The map \(\mu\mapsto\mu(\mathcal E)\) is consequently a
linear isomorphism from \(\mathcal D_w(V)\) onto \(\mathbb R\).

\begin{definition}[Bundle of densities of weight $w$ on a manifold]\label{def:variedades-riemannianas-haz-de-densidades-de-peso-w-en-una-variedad}\index{density bundle!of weight \(w\)}
Let $M$ be a smooth manifold with or without boundary. The \textbf{bundle of densities of weight $w$ on $M$}, denoted by $\boldsymbol{\mathcal{D}}_{w}(M)$, is the set
\[
\boldsymbol{\mathcal{D}}_{w}(M) := \coprod_{p \in M} \mathcal{D}_{w}(T_{p}M).
\]
$\boldsymbol{\mathcal{D}}_{w}(M)$ is a smooth vector bundle of rank $1$. Over any coordinate chart $(U, x^i)$, a local basis for this bundle is given by $|\mathbf{d}x^{1} \wedge \cdots \wedge \mathbf{d}x^{n}|^{w}$.

A \textbf{density of weight $w$ on $M$} is a smooth section of $\boldsymbol{\mathcal{D}}_{w}(M)$.
\end{definition}
 We will define the integral of a density on a smooth manifold. For this, we use a construction analogous to the one for differential forms.

 \begin{definition}\label{def:variedades-riemannianas-integral-conjunto-jordan-medible-densidad}\index{integral of densities on Euclidean domains}
 If $D\subseteq\mathbb{H}^{n}$ is a Jordan-measurable set and $\boldsymbol{\mu}$ is a continuous density of weight one on $\overline{D}$, we may write $\boldsymbol{\mu}=f|\mathbf{d}x^{1}\wedge \cdots\wedge \mathbf{d}x^{n}|$ for a unique continuous function $f\colon \overline{D}\longrightarrow \mathbb{R}$. Define the integral of $\boldsymbol{\mu}$ over $D$ by \[\int_{D}\boldsymbol{\mu}=\int_{D}f\,dx^{1}\cdots dx^{n}.\] Similarly, if $U\subseteq \mathbb{H}^{n}$ is open and $\boldsymbol{\mu}$ is a continuous density with compact support contained in $U$, define \[\int_{U}\boldsymbol{\mu}=\int_{D}\boldsymbol{\mu}\] where $D$ is a Jordan-measurable set such that
\[
\supp(\boldsymbol{\mu})\subseteq D
\quad\text{and}\quad
\overline D\subseteq U.
\]
Such a set exists: cover the support by finitely many sufficiently small closed
rectangular boxes contained in \(U\)
and take their union.
 \end{definition}
 \begin{definition}\label{def:variedades-riemannianas-integral-de-sobre-2}\index{integral of densities on manifolds}
 Let $M$ be a smooth manifold with or without boundary (not necessarily orientable). If $\boldsymbol{\mu}$ is a compactly supported density on $M$ whose support is contained in the domain of a smooth chart $(U,\phi)$, define the integral of $\boldsymbol{\mu}$ over $M$ to be \[\int_{M}\boldsymbol{\mu}=\int_{\phi(U)}(\phi^{-1})^{*}\boldsymbol{\mu}.\] We may extend this definition to arbitrary compactly supported densities by defining the \textbf{integral of $\boldsymbol{\mu}$ over $M$} by \[\int_{M}\boldsymbol{\mu}=\displaystyle\sum_{i=1}^{k}\int_{M}\psi_{i}\boldsymbol{\mu},\] where the chart domains \(U_1,\dots,U_k\) cover
\(\supp(\boldsymbol\mu)\), each
\(\supp(\psi_i)\) is contained in \(U_i\), and
\(\displaystyle \sum_{i=1}^k\psi_i=1\) on a neighborhood of
\(\supp(\boldsymbol\mu)\). These functions are obtained by applying a
partition of unity to the cover
\(U_1,\dots,U_k,M\setminus\supp(\boldsymbol\mu)\).
 \end{definition}

\begin{proposition}[Well-definedness of the integral of densities]
\label{prop:integral-densidades-bien-definida}
The preceding integral is independent of the chosen Jordan-measurable set, chart,
finite cover, and partition of unity.
\end{proposition}
\begin{proof}
Begin with an open set \(U\subseteq\mathbb H^n\) and a continuous density
\(\boldsymbol\mu=f\,|\mathbf d x^1\wedge\cdots\wedge\mathbf d x^n|\)
whose support \(K\) is compact and contained in \(U\).
The extension of \(f\) by zero outside \(U\) is continuous on a neighborhood
of any Jordan-measurable set \(D\) such that \(K\subseteq D\) and
\(\overline D\subseteq U\). If \(D_1\) and \(D_2\) are two such sets, the extended
function vanishes on both \((D_1\cup D_2)\setminus D_1\) and
\((D_1\cup D_2)\setminus D_2\). Since \(D_1\cup D_2\) is
Jordan-measurable, additivity of the Riemann integral gives
\[
\int_{D_1}f\,d\lambda_n
=\int_{D_1\cup D_2}f\,d\lambda_n
=\int_{D_2}f\,d\lambda_n.
\]
This proves independence of \(D\).

Now let \((U,\phi)\) and \((V,\psi)\) be two charts whose domains contain
\(\operatorname{supp}(\boldsymbol\mu)\), and set
\(\Theta:=\psi\circ\phi^{-1}\) on \(\phi(U\cap V)\). If
\[
(\phi^{-1})^*\boldsymbol\mu
=f_\phi\,|\mathbf d x^1\wedge\cdots\wedge\mathbf d x^n|,
\qquad
(\psi^{-1})^*\boldsymbol\mu
=f_\psi\,|\mathbf d y^1\wedge\cdots\wedge\mathbf d y^n|,
\]
the transformation law for densities gives
\[
f_\phi(x)
=f_\psi(\Theta(x))\,|\det D\Theta(x)|.
\]
The change-of-variables theorem
\ref{teorema de cambio de variable}, applied on \(\phi(U\cap V)\),
then implies
\[
\int_{\phi(U)}f_\phi\,d\lambda_n
=\int_{\psi(V)}f_\psi\,d\lambda_n,
\]
since both integrands vanish outside the image of the support. Thus the
local integral is independent of the chart.

Finally, let \((\psi_i)_{i=1}^k\) and
\((\eta_j)_{j=1}^{\ell}\) be two partitions of unity subordinate to
respective finite covers of the support of
\(\boldsymbol\mu\) by chart domains. Additivity of the local integral and the chart
independence already proved allow us to integrate each
\(\psi_i\eta_j\boldsymbol\mu\) in either of the two charts
containing its support. Since the sums are finite,
\[
\begin{aligned}
\displaystyle\sum_{i=1}^k\int_M\psi_i\boldsymbol\mu
&=\displaystyle\sum_{i=1}^k\displaystyle\sum_{j=1}^{\ell}
  \int_M\psi_i\eta_j\boldsymbol\mu\\
&=\displaystyle\sum_{j=1}^{\ell}\int_M\eta_j\boldsymbol\mu.
\end{aligned}
\]
This proves independence of both the partition and the
cover.
\end{proof}

 \begin{proposition}[The Riemannian density]\label{prop:variedades-riemannianas-la-densidad-riemanniana}\index{density!Riemannian}
 If $(M,\mathbf{g})$ is a Riemannian manifold with or without boundary, then there exists a unique smooth positive density $\boldsymbol{\mu}_{\mathbf{g}}$ defined on $M$, called the \textbf{Riemannian density}, with the property that $\boldsymbol{\mu}_{\mathbf{g}}(\mathbf{E}_{1},\dots,\mathbf{E}_{n})=1$ for every local orthonormal frame $(\mathbf{E}_{i})$.
 \end{proposition}
\begin{proof}
In a chart $(U,x^1,\dots,x^n)$, set
\[
 \boldsymbol{\mu}_{\mathbf{g}}=\sqrt{\det(\mathbf{g})}\,
 |\mathbf{d}x^1\wedge\cdots\wedge \mathbf{d}x^n|.
\]
We have the identity
\[
 \det \widetilde G
 =\det\left(\frac{\partial x}{\partial y}\right)^2\det G,
\]
where $G=(g_{ij})$ and $\widetilde G=(\widetilde g_{ij})$ are the matrices of
components in the coordinates $x$ and $y$, respectively. This identity,
together with the transformation law for a density, shows that these expressions
agree on overlaps, regardless of the sign of the Jacobian. Therefore,
they define a global smooth positive density. If $(\mathbf{E}_1,\dots,\mathbf{E}_n)$ is an
orthonormal frame and $B$ is its matrix of components in the coordinate frame,
then $B^T(g_{ij})B=I$ and
$|\det B|=(\det(\mathbf{g}))^{-\frac12}$; this gives
$\boldsymbol{\mu}_{\mathbf{g}}(\mathbf{E}_1,\dots,\mathbf{E}_n)=1$. Every density with this property agrees with
$\boldsymbol{\mu}_{\mathbf{g}}$ on each frame and consequently equals it.
\end{proof}

 We define the integral of a compactly supported function on a Riemannian manifold, not necessarily orientable, as follows:
 \begin{definition}\label{def:variedades-riemannianas-integral-de-sobre-3}\index{integral with respect to the Riemann Lebesgue measure@integral with respect to the Riemann--Lebesgue measure} Let $(M,\mathbf{g})$ be a Riemannian manifold with or without boundary.
 If $f\colon M\longrightarrow \mathbb{R}$ is a continuous compactly supported function, define the \textbf{integral of $f$ over $M$} to be \[\int_{M}f=\int_{M}f\boldsymbol{\mu}_{\mathbf{g}}.\]
 \end{definition}

Stokes's theorem turns an exterior derivative into an integral
over the boundary. Its local proof is the fundamental theorem of
calculus; a partition of unity extends it to the whole manifold.
In dimension zero, integration on an oriented manifold means the
sum of the function values multiplied by the orientation signs of
the points; a compact support contains only finitely many points.

\begin{theorem}[Stokes]
\label{teo:variedades-riemannianas-de-stokes}
\index{Stokes theorem}
Let $M$ be an oriented smooth manifold of dimension $n\geq1$,
with or without boundary, and let $\boldsymbol\omega$ be a smooth
compactly supported $(n-1)$-form. Give $\partial M$ the induced
orientation with the outward-normal-first convention, and denote
its inclusion by $\iota_\partial$. Then
\[
\int_M d\boldsymbol\omega
=\int_{\partial M}\iota_\partial^*\boldsymbol\omega.
\]
If $\partial M=\varnothing$, the right-hand side is zero.
\end{theorem}
\begin{proof}
First consider a form whose compact support lies in a single chart.
Extend its coefficients by zero outside the coordinate domain,
within $\mathbb R^n$ or $\mathbb H^n$. This extension is smooth
up to the boundary because the support does not reach the relative
edge of the chart. In coordinates, write
\[
\boldsymbol\omega=\displaystyle\sum_{i=1}^n(-1)^{i-1}a_i\,
\mathbf d x^1\wedge\cdots\wedge\widehat{\mathbf d x^i}
\wedge\cdots\wedge\mathbf d x^n.
\]
On differentiation, only the derivative of $a_i$ with respect to
$x^i$ survives, and reordering the differentials cancels $(-1)^{i-1}$:
\[
d\boldsymbol\omega=
\left(\displaystyle\sum_{i=1}^n\frac{\partial a_i}{\partial x^i}\right)
\mathbf d x^1\wedge\cdots\wedge\mathbf d x^n.
\]
On $\mathbb R^n$, each integral of $\partial a_i/\partial x^i$
vanishes by iterated integration and the fundamental theorem of
calculus, since $a_i$ has compact support. There is no boundary term.

On $\mathbb H^n=\mathbb R^{n-1}\times[0,\infty)$, the tangential
derivatives $i<n$ integrate to zero for the same reason. For the
normal derivative, with $x'=(x^1,\ldots,x^{n-1})$,
\[
\int_{\mathbb R^{n-1}}\int_0^{\infty}
\frac{\partial a_n}{\partial x^n}(x',x^n)\,dx^n\,dx'
=-\int_{\mathbb R^{n-1}}a_n(x',0)\,dx'.
\]
The outward coordinate normal is $-\boldsymbol\partial_n$.
With the standard orientation on the half-space, the induced
orientation on $x^n=0$ is $(-1)^n$ times the standard orientation
of $(x^1,\ldots,x^{n-1})$: this is the sign of the basis
$(-\boldsymbol\partial_n,\boldsymbol\partial_1,\ldots,
\boldsymbol\partial_{n-1})$ relative to the ambient basis.
The pullback of the form to the boundary is
\[
(-1)^{n-1}a_n(x',0)\,
\mathbf d x^1\wedge\cdots\wedge\mathbf d x^{n-1}.
\]
Its integral in the induced orientation is therefore the same
negative integral obtained above. When $n=1$, empty products and
integration on $\mathbb R^0$ mean evaluation; the formula gives
the negative sign of the initial endpoint. If the chart reverses
orientation, both integrals acquire the same sign, so the local
identity remains valid.

In general, choose a finite partition near the support as in the
definition of integration. Write
$\boldsymbol\omega=\displaystyle\sum_{i=1}^k\psi_i\boldsymbol\omega$.
Each summand has compact support in a chart, so the local case
applies. Sum these identities and use
\[
d\boldsymbol\omega=\displaystyle\sum_{i=1}^k d(\psi_i\boldsymbol\omega),
\qquad
\iota_\partial^*\boldsymbol\omega
=\displaystyle\sum_{i=1}^k\iota_\partial^*(\psi_i\boldsymbol\omega)
\]
to obtain the global formula. Interior charts create no artificial
boundaries: their coefficients vanish outside a compact subset of
the chart domain.
\end{proof}

\section{The Riemann--Lebesgue measure}

Integration of smooth forms suffices for elementary differential geometry, but analysis requires integrating measurable functions and working with convergence, $L^p(M)$ spaces, and weak derivatives. The Riemann--Lebesgue measure extends geometric integration to this setting while retaining, in each chart, the density determined by the metric. It allows us to handle much less regular functions while preserving the intrinsic structure of the manifold.

\begin{definition}\label{def:variedades-riemannianas-lebesgue-medible}\index{Lebesgue measurable}
 Let $M$ be a smooth manifold with or without boundary. A set $A\subseteq M$
 is \textbf{Lebesgue measurable} if, for every smooth chart
 $(U,\phi)$ of $M$, the set $\phi(A\cap U)$ is Lebesgue measurable in
 $\mathbb R^n$. Write
 \[
 \mathcal L_M:=\{A\subseteq M\mid A\text{ is Lebesgue measurable}\}.
 \]
\end{definition}
 \begin{proposition}\label{prop:variedades-riemannianas-algebra-contiene-borelianos}
 Let $M$ be a smooth manifold with or without boundary. Then
 $\mathcal{L}_{M}$ is a $\sigma$-algebra on $M$ containing the
 Borel sets.
 \end{proposition}
 \begin{proof}
 Let $(U,\phi)$ be a chart of $M$. If $A\in \mathcal{L}_{M}$, then $\phi(A\cap U)\in\mathcal{L}(n)$. Since $\mathcal{L}(n)$ is a $\sigma$-algebra, we have $\phi((M\setminus A)\cap U)=\phi(U)\setminus \phi(A\cap U)\in\mathcal{L}(n)$, and thus $M\setminus A\in \mathcal{L}_{M}$. If $(A_{k})_{k\in\mathbb{N}}\subseteq \mathcal{L}_{M}$, then $\phi(A_{k}\cap U)\in\mathcal{L}(n)$, so \[\phi\left(\left(\bigcup_{k\in\mathbb{N}}A_{k}\right)\cap U\right)=\phi\left(\bigcup_{k\in\mathbb{N}}(A_{k}\cap U)\right)=\bigcup_{k\in\mathbb{N}}\phi(A_{k}\cap U)\in \mathcal{L}(n), \] which implies $\displaystyle\bigcup_{k\in\mathbb{N}}A_{k}\in\mathcal{L}_{M}$. Now $\phi(M\cap U)=\phi(U)\in \mathcal{L}(n)$ because it is open, so $M\in \mathcal{L}_{M}$, and we conclude that $\mathcal{L}_{M}$ is a $\sigma$-algebra on $M$.

 Finally, if $O\subseteq M$ is open, then $\phi(O\cap U)\in\mathcal{L}(n)$ because it is open in $\mathbb{H}^{n}$. Therefore, $\mathcal{B}(M)\subseteq \mathcal{L}_{M}$.
 \end{proof}
 \begin{definition}\label{def:variedades-riemannianas-variedad-6}\index{local Riemann Lebesgue measure@local Riemann--Lebesgue measure} Let $(M,\mathbf{g})$ be a Riemannian manifold with or without boundary and let $(U,\phi)$ be a chart of $M$.
 If $A\in\mathcal{L}_{M}$ and $A\subseteq U$, define $\lambda_{\mathbf{g},U}(A):=\displaystyle\int_{\phi(A)}\phi_{*}\left(\sqrt{\det(\mathbf{g})}\right)d\lambda_{n}$, where $\phi_{*}\left(\sqrt{\det(\mathbf{g})}\right)=\left(\sqrt{\det(\mathbf{g})}\right)\circ \phi^{-1}$.
 \end{definition}
 \begin{proposition}\label{medida de lebesgue riemann no depende de la carta}

 \begin{enumerate}[label=(\alph*)]
 Let $(M,\mathbf{g})$ be a Riemannian manifold with or without boundary and let $A\in \mathcal{L}_{M}$.
 \item The definition of $\lambda_{\mathbf{g},U}$ is independent of the chart $(U,\phi)$ when $A\subseteq U$.
 \item $\lambda_{\mathbf{g},U}$ defines a measure on the measurable space $(U,(\mathcal{L}_{M})\restriction_{U})$.
 \end{enumerate}
 \end{proposition}
 \begin{proof}
Let $(U,\phi)$ be a chart and write
\[
G_\phi=(g_{ij}),\qquad
\rho_\phi:=\sqrt{\det G_\phi}\circ\phi^{-1}.
\]
The function $\rho_\phi$ is continuous and strictly positive on $\phi(U)$.

\begin{enumerate}[label=(\alph*)]
\item Let $(U,\phi)$ and $(V,\psi)$ be two charts with $A\subseteq U\cap V$, and
let
\[
f:=\psi\circ\phi^{-1}\colon \phi(U\cap V)\longrightarrow\psi(U\cap V)
\]
be the change of coordinates. If $G_\psi=(\widetilde g_{ab})$ is the matrix of
components of $\mathbf{g}$ in the chart $\psi$, the transformation law for
a covariant tensor of order two gives
\[
G_\psi(f(x))
=Df(x)^{-\mathsf T}G_\phi(x)Df(x)^{-1}.
\]
Taking determinants and square roots gives
\begin{equation}
\rho_\phi(x)
=\rho_\psi(f(x))\,|\det Df(x)|.
\label{eq:transformacion-densidad-riemann-lebesgue}
\end{equation}
The change-of-variables formula applied to the nonnegative function
$\rho_\psi\mathbf 1_{\psi(A)}$ implies
\[
\begin{aligned}
\lambda_{\mathbf{g},V}(A)
&=\int_{\psi(A)}\rho_\psi(y)\,d\lambda_n(y)\\
&=\int_{\phi(A)}\rho_\psi(f(x))|\det Df(x)|\,d\lambda_n(x)\\
&=\int_{\phi(A)}\rho_\phi(x)\,d\lambda_n(x)
 =\lambda_{\mathbf{g},U}(A),
\end{aligned}
\]
where the last equality uses
\eqref{eq:transformacion-densidad-riemann-lebesgue}. Thus the local
measure is independent of the chart.

\item For every Lebesgue measurable set
$B\subseteq\phi(U)$, define
\[
\nu_\phi(B):=\int_B\rho_\phi\,d\lambda_n.
\]
Since $\rho_\phi\geq0$ is measurable, $\nu_\phi$ is a measure:
$\sigma$-additivity follows by applying the monotone convergence theorem
to sums of indicator functions of pairwise disjoint sets.
Moreover,
\[
\lambda_{\mathbf{g},U}(A)=\nu_\phi(\phi(A)),
\qquad A\in(\mathcal L_M)\restriction_U.
\]
Since $\phi$ is a measurable bijection with measurable inverse, this
identity transfers $\sigma$-additivity from $\nu_\phi$ to
$\lambda_{\mathbf{g},U}$. Thus, $\lambda_{\mathbf{g},U}$ is a measure on
$(U,(\mathcal L_M)\restriction_U)$.
\end{enumerate}
 \end{proof}
 \begin{definition}\label{medibles ajenos dos a dos, definicion de medida Riemann--Lebesgue}\index{Riemann Lebesgue measure@Riemann--Lebesgue measure}\index{manifold!Riemann--Lebesgue measure}\glsadd{medida-riemann-lebesgue}
 Let $(M,\mathbf{g})$ be a Riemannian manifold with or without boundary. For each
 $A\in\mathcal L_M$, choose a countable atlas
 $\mathcal A=\{(U_k,\phi_k)\mid k\in\mathbb N\}$ of $M$ and define
 \[
 A_0:=A\cap U_0,
 \qquad
 A_{k+1}:=(A\cap U_{k+1})\setminus\bigcup_{j=0}^{k}A_j.
 \]
 The sets $A_k$ are measurable, pairwise disjoint, satisfy
 $A_k\subseteq U_k$, and cover $A$. Define
 \[
 \lambda_{\mathbf{g}}(A):=\displaystyle\sum_{k=0}^{\infty}
 \lambda_{\mathbf{g},U_k}(A_k).
 \]
 This function is called the \textbf{Riemann--Lebesgue measure} associated with
 $\mathbf{g}$.
 \end{definition}
\begin{figura}[h]
    \includegraphics[scale=1]{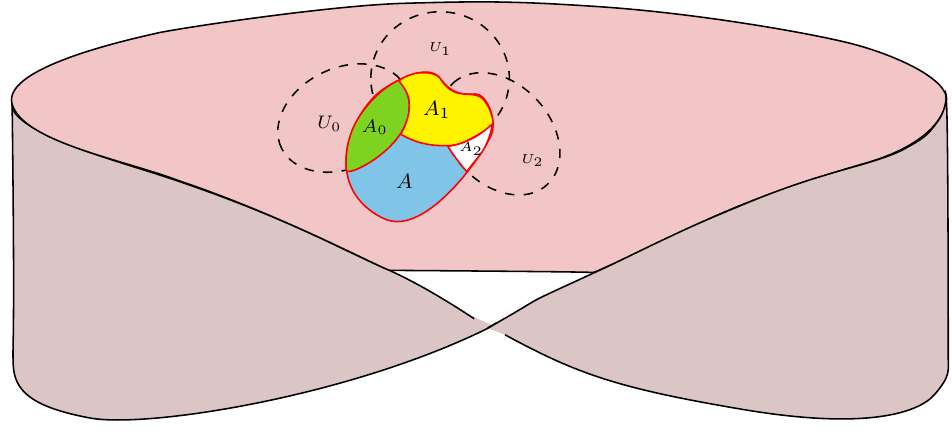}
    \caption{The definition does not depend on the orientability of the manifold.}
    \label{fig:mobius}
\end{figura}

 \begin{proposition}\label{prop:variedades-riemannianas-definicion-no-depende-eleccion-atlas-numerable}
 Let $(M,\mathbf{g})$ be a Riemannian manifold with or without boundary. The
 definition of $\lambda_{\mathbf{g}}$ is independent of the choice of countable
 atlas.
 \end{proposition}
 \begin{proof}
Let $\widetilde{\mathcal{A}}$ be another countable atlas of $M$. Let $\widetilde{A}_{0}:= A\cap \widetilde{U}_{0}$; for each $i\in\mathbb{N}$ with $i\geq 0$, define $\widetilde{A}_{i+1}:=(A\cap \widetilde{U}_{i+1})\setminus \displaystyle\bigcup_{j=0}^{i}\widetilde{A}_{j}$. By Proposition~\ref{medida de lebesgue riemann no depende de la carta}, $\lambda_{\mathbf{g},U_{k}}$ and $\lambda_{\mathbf{g},\widetilde{U}_{i}}$ are measures on $U_{k}$ and $\widetilde{U}_{i}$, respectively, so they are $\sigma$-additive, and using $\displaystyle\bigcup_{k\in\mathbb{N}}A_{k}=\displaystyle\bigcup_{i\in\mathbb{N}}\widetilde{A}_{i}=A$, \[\lambda_{\mathbf{g},U_{k}}(A_{k})=\lambda_{\mathbf{g},U_{k}}(A_{k}\cap A)=\lambda_{\mathbf{g},U_{k}}\left(A_{k}\cap \bigcup_{i\in\mathbb{N}}\widetilde{A}_{i}\right)=\lambda_{\mathbf{g},U_{k}}\left(\bigcup_{i\in\mathbb{N}}(A_{k}\cap\widetilde{A}_{i})\right)=\displaystyle\sum_{i\in\mathbb{N}}\lambda_{\mathbf{g},U_{k}}(A_{k}\cap\widetilde{A}_{i}).\]
Apply the same decomposition to $\lambda_{\mathbf{g},\widetilde{U}_{i}}(\widetilde{A}_{i})=\displaystyle\sum_{k\in\mathbb{N}}\lambda_{\mathbf{g},\widetilde{U}_{i}}(\widetilde{A}_{i}\cap A_{k})$. By Proposition~\ref{medida de lebesgue riemann no depende de la carta}, we know that $\lambda_{\mathbf{g},U_{k}}(A_{k}\cap\widetilde{A}_{i})=\lambda_{\mathbf{g},\widetilde{U}_{i}}(\widetilde{A}_{i}\cap A_{k})$ since $A_{k}\cap\widetilde{A}_{i}\subseteq U_{k}\cap\widetilde{U}_{i}$.

Therefore, \[\displaystyle\sum_{k=0}^{\infty}\lambda_{\mathbf{g},U_{k}}(A_{k})=\displaystyle\sum_{k=0}^{\infty}\displaystyle\sum_{i=0}^{\infty}\lambda_{\mathbf{g},U_{k}}(A_{k}\cap \widetilde{A}_{i})\] \[=\displaystyle\sum_{k=0}^{\infty}\displaystyle\sum_{i=0}^{\infty} \lambda_{\mathbf{g},\widetilde{U}_{i}}(\widetilde{A}_{i}\cap A_{k})=\displaystyle\sum_{i=0}^{\infty}\displaystyle\sum_{k=0}^{\infty}\lambda_{\mathbf{g},\widetilde{U}_{i}}(\widetilde{A}_{i}\cap A_{k})=\displaystyle\sum_{i=0}^{\infty}\lambda_{\mathbf{g},\widetilde{U}_{i}}(\widetilde{A}_{i}).\]
 \end{proof}
The measure $\lambda_{\mathbf{g}}$ agrees with $\lambda_{\mathbf{g},U}$ on measurable sets contained in $U$:
 \begin{proposition}\label{prop:variedades-riemannianas-variedad-riemanniana-frontera-carta-suave}
 Let $(M,\mathbf{g})$ be a Riemannian manifold with or without boundary and let $(U,\phi)$ be a smooth chart. If $A\in \mathcal{L}_{M}$ and $A\subseteq U$, then $\lambda_{\mathbf{g}}(A)=\lambda_{\mathbf{g},U}(A)$.
 \end{proposition}
\begin{proof}
 Let $(U_{k},\phi_{k})_{k\in \mathbb{N}}$ be a countable atlas of $M$ such that $U_{0}=U$. Consider $A_{0}=A\cap U_{0}=A\cap U=A$ and $A_{k+1}=(A\cap U_{k+1})\setminus \displaystyle\bigcup_{j=0}^{k}A_{j}$. We know that $(A_{k})_{k\in \mathbb{N}}$ is a sequence of pairwise disjoint measurable sets whose union is $A$. Since $A_{0}=A$, observe that this implies $A_{k}=\varnothing$ for every $k\geq 1$. Therefore, \[\lambda_{\mathbf{g}}(A)=\displaystyle\sum_{k=0}^{\infty}\lambda_{\mathbf{g},U_{k}}(A_{k})=\lambda_{\mathbf{g},U_{0}}(A_{0})=\lambda_{\mathbf{g},U}(A).\]
\end{proof}

 \begin{theorem}\label{teo: la medida de riemann lebesgue es de radon}\index{Riemann Lebesgue measure@Riemann--Lebesgue measure!Radon property}
 Let $(M,\mathbf{g})$ be a Riemannian manifold with or without boundary. Then
 $\lambda_{\mathbf{g}}$ is a complete Radon measure on
 $(M,\mathcal{L}_{M})$.
 \end{theorem}
 \begin{proof}
Fix a countable atlas $\{(U_j,\phi_j)\mid j\in\mathbb N_0\}$.
For each $B\in\mathcal L_M$, set
\[
B^{(0)}=B\cap U_0,
\qquad
B^{(j+1)}=(B\cap U_{j+1})\setminus\bigcup_{l=0}^{j}B^{(l)}.
\]
These pieces are measurable, pairwise disjoint, cover $B$, and satisfy
$B^{(j)}\subseteq U_j$. If $(A_k)_{k=1}^{\infty}$ is a countable family of
pairwise disjoint measurable sets and $\displaystyle A=\bigcup_{k=1}^{\infty}A_k$, then, for each $j$,
\[
A^{(j)}=\bigcup_{k=1}^{\infty} A_k^{(j)},
\]
and the union on the right is also disjoint. Countable additivity
of the local measures and interchange of two series of nonnegative
terms give
\[
\begin{aligned}
\lambda_{\mathbf{g}}(A)
&=\displaystyle\sum_{j=0}^{\infty}\lambda_{\mathbf{g},U_j}(A^{(j)})
=\displaystyle\sum_{j=0}^{\infty}\displaystyle\sum_{k=1}^{\infty}\lambda_{\mathbf{g},U_j}(A_k^{(j)})\\
&=\displaystyle\sum_{k=1}^{\infty}\displaystyle\sum_{j=0}^{\infty}\lambda_{\mathbf{g},U_j}(A_k^{(j)})
=\displaystyle\sum_{k=1}^{\infty}\lambda_{\mathbf{g}}(A_k).
\end{aligned}
\]
Moreover, $\lambda_{\mathbf{g}}(\varnothing)=0$; therefore, $\lambda_{\mathbf{g}}$ is a measure.

We first establish local regularity. In a chart $(U,\phi)$, the
local measure is the transport by $\phi$ of the measure
\[
\nu(B)=\int_B w\,d\lambda_n,
\qquad
w=\phi_*\sqrt{\det(\mathbf{g})},
\]
on the relatively open set $\phi(U)\subseteq\mathbb H^n$. The function $w$
is continuous and strictly positive. Take an increasing sequence of
relatively compact open sets $W_m$, $m\geq1$, whose union is $\phi(U)$, and put $W_0=\varnothing$. On each
$\overline W_m$, the function $w$ has a finite upper bound and a positive
lower bound. Regularity of Lebesgue measure, applied on
each $W_m$, then shows that for every measurable set
$B\subseteq\phi(U)$,
\[
\nu(B)=\displaystyle\sup\{\nu(K)\mid K\subseteq B, K\text{ compact}\}
=\displaystyle\inf\{\nu(O)\mid B\subseteq O, O\text{ open in }\phi(U)\}.
\]
Indeed, for inner approximation, first choose $m$ so that
$\nu(B\setminus W_m)$ is small and then a compact set
$K\subseteq B\cap W_m$ whose complement has small Lebesgue measure;
the upper bound for $w$ turns this last estimate into an
estimate for $\nu$. If $\nu(B)=\infty$, the same argument applied to
$B\cap W_m$ gives compact sets of arbitrarily large measure. For outer
approximation, decompose $B$ into the measurable sets
$B\cap(W_m\setminus W_{m-1})$ and, again using the upper bound for $w$
on $\overline W_m$, choose for the $m$th piece an open set whose excess
has measure less than $\frac{\varepsilon}{2^{m+1}}$. The union of these open sets
contains $B$ and its total excess has measure less than $\varepsilon$.

We may choose the countable atlas so that each $U_j$ is a precompact
coordinate domain with bounded image and
$\lambda_{\mathbf{g}}(U_j)<\infty$. For $A\in\mathcal L_M$, let $(A_j)$ be the disjoint
pieces associated with this atlas. If $\lambda_{\mathbf{g}}(A)<\infty$ and
$\varepsilon>0$, choose $N$ so that
\[
\displaystyle\sum_{j>N}\lambda_{\mathbf{g}}(A_j)<\frac{\varepsilon}{2}.
\]
Local regularity gives compact sets $K_j\subseteq A_j$, for
$0\leq j\leq N$, such that
\[
\displaystyle\sum_{j=0}^{N}\lambda_{\mathbf{g}}(A_j\setminus K_j)<\frac{\varepsilon}{2}.
\]
The finite union $\displaystyle K=\bigcup_{j=0}^{N}K_j$ is compact, is contained in
$A$, and satisfies $\lambda_{\mathbf{g}}(A\setminus K)<\varepsilon$. If
$\lambda_{\mathbf{g}}(A)=\infty$, first take a partial sum greater than a given
number and approximate its terms from within; this gives compact sets
contained in $A$ of arbitrarily large measure. This proves inner
regularity.

For outer regularity, we need only consider
$\lambda_{\mathbf{g}}(A)<\infty$. For each $j$, use local regularity to
choose an open subset $O_j$ of $M$ such that
\[
A_j\subseteq O_j\subseteq U_j,
\qquad
\lambda_{\mathbf{g}}(O_j\setminus A_j)<\frac{\varepsilon}{2^{j+1}}.
\]
Then $\displaystyle O=\bigcup_{j=0}^{\infty}O_j$ is open, contains $A$, and
\[
\lambda_{\mathbf{g}}(O\setminus A)
\leq\displaystyle\sum_{j=0}^{\infty}\lambda_{\mathbf{g}}(O_j\setminus A_j)<\varepsilon.
\]
If $\lambda_{\mathbf{g}}(A)=\infty$, every open set containing $A$ has infinite
measure. This proves outer regularity.

Finally, every point has a coordinate domain $V$ with compact
closure contained in another coordinate domain. The function $w$ is bounded
on the image of this closure, and $\phi(V)$ can be chosen bounded; therefore,
$\lambda_{\mathbf{g}}(V)<\infty$. The measure is locally finite. It is also complete:
if $N\subseteq Z$ and $\lambda_{\mathbf{g}}(Z)=0$, positivity of $w$ and its lower
bounds on compact sets imply that $\phi(Z\cap U)$ has Lebesgue measure
zero in every chart. Every subset of this set is
Lebesgue measurable and has local measure zero, so
$N\in\mathcal L_M$ and $\lambda_{\mathbf{g}}(N)=0$. We conclude that $\lambda_{\mathbf{g}}$ is a
Radon measure on the indicated measurable space.
 \end{proof}
\begin{remark}\label{obs:notacion-integral-medida-riemanniana}
The expression $\displaystyle \int_M f\,d\lambda_{\mathbf{g}}$ denotes the Lebesgue integral with respect
to the measure $\lambda_{\mathbf{g}}$. For every continuous compactly supported function,
\[
 \int_M f\,d\lambda_{\mathbf{g}}=\int_M f\,\boldsymbol{\mu}_{\mathbf{g}}.
\]
If $M$ is oriented and $dV_{\mathbf{g}}$ corresponds to the chosen orientation,
$\boldsymbol{\mu}_{\mathbf{g}}$ is the density associated with $dV_{\mathbf{g}}$, and therefore
\[
 \int_M f\,d\lambda_{\mathbf{g}}=\int_M f\,dV_{\mathbf{g}}.
\]
Thus, $dV_{\mathbf{g}}$, $\boldsymbol{\mu}_{\mathbf{g}}$, and $d\lambda_{\mathbf{g}}$ denote, respectively, a form,
a density, and the differential notation for a measure.
\end{remark}
 Local finiteness of the Riemann--Lebesgue measure has the following consequence:
 \begin{remark}\label{M compacta tiene medida finita}
 Let $(M,\mathbf{g})$ be a compact Riemannian manifold with or without boundary. Then $\lambda_{\mathbf{g}}(M)<\infty$.
 \end{remark}
 \begin{proof}
 Since $M$ is compact and $\lambda_{\mathbf{g}}$ is locally finite, we can find a finite open cover of $M$, say $M=\displaystyle\bigcup_{j=0}^{N}U_{j}$, such that $\lambda_{\mathbf{g}}(U_{j})<\infty$ for every $j\in\{0,\dots,N\}$.

 Then \[\lambda_{\mathbf{g}}(M)\leq \displaystyle\sum_{j=0}^{N}\lambda_{\mathbf{g}}(U_{j})<\infty.\]
 \end{proof}
 For later analytic applications, we will explicitly need
the Borel completion and $\sigma$-finiteness of the Riemann--Lebesgue measure.

\begin{lemma}[Borel completion and $\sigma$-finiteness of the Riemann--Lebesgue measure]
\label{lem:medida-riemanniana-completacion-borel-sigma-finita}
\index{Riemann Lebesgue measure@Riemann--Lebesgue measure!Borel completion and sigma-finiteness}
Let $(M,\mathbf{g})$ be a Riemannian manifold with or without boundary and let
$\lambda_{\mathbf{g}}$ be the associated Riemann--Lebesgue measure. Then
\[
\mathcal L_M=\overline{\mathcal B(M)}^{\,\lambda_{\mathbf{g}}},
\]
where the right-hand side denotes the completion of the Borel $\sigma$-algebra
with respect to $\lambda_{\mathbf{g}}$. Moreover,
$(M,\mathcal L_M,\lambda_{\mathbf{g}})$ is $\sigma$-finite.
\end{lemma}

\begin{proof}
By Proposition~\ref{bolacoordenadaregular}, we may choose a countable
cover $(U_j,\phi_j)_{j\in\mathbb N}$ by regular coordinate balls and,
in the boundary case, by the corresponding regular coordinate
half-balls, so that $\overline{U_j}$ is compact for every $j$.

We first prove equality of the $\sigma$-algebras. The inclusion
$\overline{\mathcal B(M)}^{\,\lambda_{\mathbf g}}\subseteq\mathcal L_M$
follows from Theorem~\ref{teo: la medida de riemann lebesgue es de radon}:
$\mathcal B(M)\subseteq\mathcal L_M$ and $\lambda_{\mathbf g}$ is complete on
$\mathcal L_M$.

For the reverse inclusion, take $A\in\mathcal L_M$ and form the
disjoint decomposition
\[
A_0:=A\cap U_0,
\qquad
A_j:=(A\cap U_j)\setminus\bigcup_{i=0}^{j-1}A_i,
\qquad j\geq1.
\]
Then $A=\displaystyle\bigcup_{j=0}^{\infty}A_j$, each $A_j$ belongs to
$\mathcal L_M$, and $A_j\subseteq U_j$. By the definition of $\mathcal L_M$,
$\phi_j(A_j)$ is Lebesgue measurable in $\phi_j(U_j)$. Since the Lebesgue $\sigma$-algebra
of a Euclidean open set or a relatively open subset of the
half-space is the completion of its Borel $\sigma$-algebra, there exist a
Borel set $B_j\subseteq\phi_j(U_j)$ and a set of Lebesgue measure
zero $N_j\subseteq\phi_j(U_j)$ such that
\[
\phi_j(A_j)\,\triangle\,B_j\subseteq N_j.
\]

Since $\overline{U_j}$ is compact, the coordinate density
$\sqrt{\det(\mathbf g)}\circ\phi_j^{-1}$ is bounded above and,
on the compact sets exhausting $\phi_j(U_j)$, is also bounded away from
zero. In particular, every set of Lebesgue measure zero in
$\phi_j(U_j)$ has preimage of $\lambda_{\mathbf g}$-measure zero. Therefore,
\[
Z_j:=\phi_j^{-1}(N_j)
\]
is $\lambda_{\mathbf g}$-null. Moreover,
$C_j:=\phi_j^{-1}(B_j)$ is Borel in $U_j$ and, since $U_j$ is open, it is
Borel in $M$. Defining
\[
C:=\bigcup_{j=0}^{\infty}C_j,
\qquad
Z:=\bigcup_{j=0}^{\infty}Z_j,
\]
gives $C\in\mathcal B(M)$, $\lambda_{\mathbf g}(Z)=0$, and
\[
A\triangle C\subseteq Z.
\]
Thus, $A$ belongs to the completion of $\mathcal B(M)$, proving
\[
\mathcal L_M\subseteq\overline{\mathcal B(M)}^{\,\lambda_{\mathbf g}}.
\]

It remains to prove $\sigma$-finiteness. For each $j$, compactness of
$\overline{U_j}$ and local finiteness of $\lambda_{\mathbf g}$ imply
\[
\lambda_{\mathbf g}(U_j)
\leq
\lambda_{\mathbf g}(\overline{U_j})<\infty.
\]
Since $M=\displaystyle\bigcup_{j=0}^{\infty}U_j$, the measure
$\lambda_{\mathbf g}$ is $\sigma$-finite.
\end{proof}

 The remaining properties of the Riemann--Lebesgue measure are left as an exercise:

\begin{exercise}[Basic properties of the Riemann--Lebesgue measure]\label{ej:propiedades-medida-riemann-lebesgue}
Let $(M,\mathbf{g})$ be a Riemannian manifold with or without boundary and let $\lambda_{\mathbf{g}}$ be the associated Riemann--Lebesgue measure.

\begin{enumerate}
\item Prove that $\lambda_{\mathbf{g}}$ is complete, that is, if $A\subset M$ is measurable and $\lambda_{\mathbf{g}}(A)=0$, then every subset of $A$ is measurable.

\item Let $(U,\phi)$ be a coordinate chart. Prove that if $A\subset U$ satisfies $\lambda_{\mathbf{g}}(A)=0$, then $\lambda_{n}\bigl(\phi(A)\bigr)=0$ in $\mathbb R^n$, where $\lambda_{n}$ denotes Lebesgue measure on $\mathbb{R}^{n}$.

\item Prove that a function $f\colon M\longrightarrow\mathbb R$ is $\mathcal L_M$-measurable
if and only if, for every chart $(U,\phi)$, the function
\[
f\circ \phi^{-1}\colon \phi(U)\longrightarrow\mathbb R
\]
is $\lambda_{n}$-measurable on $\phi(U)\subset\mathbb R^n$.

\item Prove that if $U\subset M$ is nonempty and open, then $\lambda_{\mathbf{g}}(U)>0$.

\item Let $(M,\mathbf{g})$ be orientable and let $dV_{\mathbf{g}}$ be the Riemannian volume form associated with an orientation. Prove that, for every continuous compactly supported function,
\[
\int_M f\,d\lambda_{\mathbf{g}}=\int_M f\,dV_{\mathbf{g}}.
\]

\end{enumerate}

\end{exercise}

\section{Gradient, divergence, and the Laplace--Beltrami operator}\label{sec:operadores-diferenciales-riemannianos}
\label{sec:operadores-grad-div-lap}

The gradient, divergence, and Laplace--Beltrami operator form the basic differential vocabulary of Riemannian analysis. They generalize familiar Euclidean operations, but their intrinsic definitions show exactly where the metric, connection, and measure enter. In Chapter~\ref{cap:operadores-diferenciales-parciales-haces}, after defining the general notion of a differential operator on bundles, we will verify that these constructions are examples of it. These operations will be used to formulate the elliptic equations, energy identities, and geometric flows of later chapters.

We first define the gradient of a scalar function $f\colon M\longrightarrow \mathbb{R}$.

The covariant derivative of $f$ equals the covector $df$. To obtain a
vector field, we must raise this index with the metric. In the Euclidean
case, this identification is expressed by
\[
\langle \operatorname{grad}f,\mathbf{X}\rangle_{\overline{\mathbf{g}}}
=\mathbf{X}(f),
\]
and leads to the following intrinsic definition.
\begin{definition}\label{gradiente}\index{gradient}
Let $(M,\mathbf{g})$ be a Riemannian manifold (with or without boundary).
For a function $f\colon M\longrightarrow\mathbb{R}$ of class $C^{1}$, the \emph{gradient} of $f$ is the unique continuous vector field $\operatorname{grad}(f)\in\Gamma^0(TM)$ such that
\[
\langle \operatorname{grad}(f),\mathbf{X}\rangle_{\mathbf{g}}=df(\mathbf{X})\quad\text{for every }\mathbf{X}\in\mathfrak{X}(M).
\]
That is, $\operatorname{grad}(f)=(df)^{\sharp}$, where $^{\sharp}\colon T^{*}M\longrightarrow TM$ is the musical isomorphism induced by $\mathbf{g}$ given in Definition~\ref{isomorfismos musicales}.
\end{definition}

\begin{remark}[Local expression for the gradient]\label{obs:variedades-riemannianas-expresion-local-del-gradiente}
Let $(U,\phi)$ be a chart with coordinates $(x^{1},\dots,x^{n})$.
In this chart, the differential of $f$ is written as
\[
df=\partial_{j}f\,\mathbf{d}x^{j},\qquad
\partial_{j}f=\displaystyle\frac{\partial(f\circ\phi^{-1})}{\partial x^{j}}.
\]
Applying the musical isomorphism $^{\sharp}\colon T^{*}M\longrightarrow TM$ gives
\[
(df)^{\sharp}=g^{ij}\partial_{j}f\,\boldsymbol{\partial}_{i}.
\]

\end{remark}

The gradient satisfies the Leibniz rule, just as the differential of a function does:
\begin{proposition}[Leibniz rule for the gradient]\label{leibniz gradiente}\index{Leibniz rule for the gradient}
 Let $(M,\mathbf{g})$ be a Riemannian manifold with or without boundary. If $f,g\in C^{1}(M)$, then \[\operatorname{grad}(fg)=f\operatorname{grad}(g)+g\operatorname{grad}(f).\]
\end{proposition}
\begin{proof}
For every smooth vector field \(\mathbf X\),
\[
\begin{aligned}
\langle\operatorname{grad}(fg),\mathbf X\rangle
&=d(fg)(\mathbf X)
=f\,dg(\mathbf X)+g\,df(\mathbf X)\\
&=\langle f\operatorname{grad}g+g\operatorname{grad}f,\mathbf X\rangle.
\end{aligned}
\]
Nondegeneracy of the metric gives equality of the fields.
\end{proof}
Divergence of a vector field is a central notion in vector calculus. We introduce its analogue on Riemannian manifolds:

\begin{definition}\label{def:variedades-riemannianas-divergencia}\index{divergence}
Let $(M,\mathbf{g})$ be a Riemannian manifold (with or without boundary) and let $\mathbf{X}$ be a vector field of class $C^{1}$.
Define the \emph{divergence} of $\mathbf{X}$ by
\[
\operatorname{div}(\mathbf{X})=\operatorname{tr}(\nabla \mathbf{X}),
\]
where $\operatorname{tr}$ denotes the natural contraction of a tensor of type $(1,1)$ in its two indices.
In a local orthonormal frame $(\mathbf{e}_1,\dots,\mathbf{e}_n)$, we recover the classical expression
\[
\operatorname{div}(\mathbf{X})=\displaystyle\sum_{i=1}^{n} \mathbf{g}(\nabla_{\mathbf{e}_i}\mathbf{X},\mathbf{e}_i).
\]
\end{definition}

In local coordinates, we have the following:

\begin{proposition}\label{prop:variedades-riemannianas-campo-vectorial-coordenadas-locales}
Let $(M,\mathbf{g})$ be a Riemannian manifold with or without boundary and let
$\mathbf{X}=X^{i}\boldsymbol{\partial}_i$ be a vector field of class $C^1$.
In local coordinates $(x^1,\dots,x^n)$, we have
\[
\operatorname{div}(\mathbf{X})=\displaystyle\frac{1}{\sqrt{\det(\mathbf{g})}}
\displaystyle\boldsymbol{\partial}_{i}\left(\sqrt{\det(\mathbf{g})} X^{i}\right).
\]
\end{proposition}

\begin{proof}
By definition, divergence is the trace of $\nabla \mathbf{X}$:
\[
\operatorname{div}(\mathbf{X})=\operatorname{tr}(\nabla \mathbf{X}).
\]
In local coordinates $(x^1,\dots,x^n)$ with coordinate frame
$\boldsymbol{\partial}_j=\boldsymbol{\partial}_{j}$, we have
\[
\nabla_{\boldsymbol{\partial}_j}\mathbf{X}
=(\partial_j X^i+\Gamma^i_{jk}X^k)\boldsymbol{\partial}_i,
\quad\text{that is}\quad
(\nabla_j \mathbf{X})^i=\partial_j X^i+\Gamma^i_{jk}X^k.
\]

Taking the trace of $\nabla \mathbf{X}$ means contracting the covariant differentiation index $j$ with the contravariant index $i$, which in index notation is written as
\[
\operatorname{div}(\mathbf{X})=\nabla_i X^i.
\]
Substituting the expression for $(\nabla_j \mathbf{X})^i$ gives
\[
\operatorname{div}(\mathbf{X})=\partial_i X^i+\Gamma^i_{ik}X^k.
\]

We now use the identity
\[
\Gamma^i_{ik}=\displaystyle\frac{1}{\sqrt{\det(\mathbf{g})}}\partial_k\big(\sqrt{\det(\mathbf{g})}\big),
\]
obtained in Remark~\ref{derivada raiz del determinante de la metrica}.

Substituting into the expression for $\operatorname{div}(\mathbf{X})$ gives
\[
\operatorname{div}(\mathbf{X})=\partial_i X^i+\displaystyle\frac{1}{\sqrt{\det(\mathbf{g})}}\partial_k\big(\sqrt{\det(\mathbf{g})}\big)X^k.
\]

Applying the Leibniz rule to the product $\sqrt{\det(\mathbf{g})}X^i$,
\[
\partial_i\big(\sqrt{\det(\mathbf{g})}X^i\big)
=\sqrt{\det(\mathbf{g})}\partial_i X^i+X^i\partial_i\sqrt{\det(\mathbf{g})},
\]
and dividing by $\sqrt{\det(\mathbf{g})}$ gives
\[
\displaystyle\frac{1}{\sqrt{\det(\mathbf{g})}}\partial_i\big(\sqrt{\det(\mathbf{g})}X^i\big)
=\partial_i X^i+\displaystyle\frac{1}{\sqrt{\det(\mathbf{g})}}\partial_i\big(\sqrt{\det(\mathbf{g})}\big)X^i.
\]

Comparing with the preceding expression yields the compact formula
\[
\operatorname{div}(\mathbf{X})=\displaystyle\frac{1}{\sqrt{\det(\mathbf{g})}}
\partial_i\big(\sqrt{\det(\mathbf{g})}X^i\big),
\]
as desired.
\end{proof}
Divergence satisfies the following Leibniz rule:
\begin{proposition}[Leibniz rule for divergence]\label{leibniz divergencia}\index{Leibniz rule for divergence}
 Let $(M,\mathbf{g})$ be a Riemannian manifold with or without boundary. Then, for every $\mathbf{X}\in \mathfrak{X}(M)$ and every $f\in C^{\infty}(M)$, \[\operatorname{div}(f\mathbf{X})=f\operatorname{div}(\mathbf{X})+\langle \operatorname{grad}(f),\mathbf{X}\rangle_{\mathbf{g}}.\]
\end{proposition}
\begin{proof}
In an orthonormal frame \((\mathbf e_i)\),
\[
\begin{aligned}
\operatorname{div}(f\mathbf X)
&=\displaystyle\sum_{i=1}^{n}\langle\nabla_{\mathbf e_i}(f\mathbf X),\mathbf e_i\rangle\\
&=\displaystyle\sum_{i=1}^{n}\mathbf e_i(f)\langle\mathbf X,\mathbf e_i\rangle
+f\displaystyle\sum_{i=1}^{n}\langle\nabla_{\mathbf e_i}\mathbf X,\mathbf e_i\rangle\\
&=\langle\operatorname{grad}f,\mathbf X\rangle
+f\operatorname{div}\mathbf X.
\end{aligned}
\]
\end{proof}
The definition extends naturally to tensor fields of type $(k,l)$:

\begin{definition}\label{def:variedades-riemannianas-variedad-7}\index{divergence of a tensor field}
Let $(M,\mathbf{g})$ be a Riemannian manifold (with or without boundary) and let $\mathbf{F}\in \Gamma(T^{(k,l)}(TM))$ with $k\geq 1$. Define
\[
\operatorname{div}(\mathbf{F})\in \Gamma(T^{(k-1,l)}(TM))
\]
to be the contraction of the differentiation index of $\nabla \mathbf{F}$ with the first contravariant index of $\mathbf{F}$.
In index notation,
\[
(\operatorname{div}\mathbf{F})^{i_2\dots i_k}_{j_1\dots j_l}
= \nabla_a F^{a i_2\dots i_k}_{j_1\dots j_l}.
\]
\end{definition}

\begin{proposition}
\label{prop:variedades-riemannianas-coordenadas-locales-su-divergencia-esta}
Let $(M,\mathbf g)$ be an $n$-dimensional Riemannian manifold with
or without boundary, and let $\mathbf F$ be a tensor of type $(k,l)$,
with $k\geq1$. In local coordinates $(x^1,\ldots,x^n)$, its divergence is
\[
\begin{aligned}
(\operatorname{div}\mathbf F)^{i_2\cdots i_k}_{j_1\cdots j_l}
={}&\displaystyle\sum_{a=1}^n
\partial_a F^{a i_2\cdots i_k}_{j_1\cdots j_l}
+\displaystyle\sum_{a,m=1}^n
\Gamma^a_{am}F^{m i_2\cdots i_k}_{j_1\cdots j_l}\\
&+\displaystyle\sum_{r=2}^k\displaystyle\sum_{a,m=1}^n
\Gamma^{i_r}_{am}F^{a i_2\cdots m\cdots i_k}_{j_1\cdots j_l}\\
&-\displaystyle\sum_{s=1}^l\displaystyle\sum_{a,m=1}^n
\Gamma^m_{a j_s}F^{a i_2\cdots i_k}_{j_1\cdots m\cdots j_l}.
\end{aligned}
\]
The free indices $i_2,\ldots,i_k,j_1,\ldots,j_l$ range over
$\{1,\ldots,n\}$. In the sum indexed by $r$, $m$ occupies
contravariant position $r$; in the sum indexed by $s$, it occupies
covariant position $s$. A sum over an empty index set is zero.
\end{proposition}
\begin{proof}
By definition, contract the derivative index with the first
contravariant index:
\[
(\operatorname{div}\mathbf F)^{i_2\cdots i_k}_{j_1\cdots j_l}
=\displaystyle\sum_{a=1}^n
(\nabla_a\mathbf F)^{a i_2\cdots i_k}_{j_1\cdots j_l}.
\]
Apply Corollary~\ref{cor:derivada-covariante-total} before summing
over $a$. The ordinary derivative produces the first term.
The connection applied to the first contravariant factor produces
\[
\displaystyle\sum_{a,m=1}^n
\Gamma^a_{am}F^{m i_2\cdots i_k}_{j_1\cdots j_l};
\]
here the upper index of the Christoffel symbol is already $a$,
since that is the component being contracted. Applying the
connection to the remaining contravariant factors, in positions
$r=2,\ldots,k$, gives the second line of the formula. Applying it
to each covariant factor, in position $s=1,\ldots,l$, gives the
third line with a minus sign from the dual connection. The Leibniz
rule for tensor products adds exactly these contributions.
\end{proof}

The Laplacian is defined by the property that characterizes it in
Euclidean space: it is the divergence of the gradient of a function.
\begin{definition}\label{def:variedades-riemannianas-operador-variedad-riemanniana-funcion-clase}\index{Laplace Beltrami operator@Laplace--Beltrami operator}
Let $(M,\mathbf{g})$ be a Riemannian manifold with or without boundary and let $f\colon M\longrightarrow \mathbb{R}$ be a function of class $C^{2}$.
The \textit{Laplace--Beltrami operator} is defined by
\[
\Delta_{\mathbf{g}} f:=\operatorname{div}(\operatorname{grad}f).
\]
\end{definition}

\begin{remark}\label{obs:variedades-riemannianas-recordemos-primero-gradiente-campo-vectorial-caracterizado}
First recall that the gradient of $f$ is the vector field characterized by
\[
\mathbf{g}(\operatorname{grad}f,\mathbf{X})=df(\mathbf{X}),\qquad \mathbf{X}\in\mathfrak{X}(M).
\]
In local coordinates $(x^1,\dots,x^n)$, this becomes
\[
\operatorname{grad}f=g^{ij}\partial_j f\,\boldsymbol{\partial}_i,
\]
where $(g^{ij})$ are the entries of the inverse metric matrix.

Applying the definition of divergence now gives, in coordinates,
\[
\operatorname{div}(\mathbf{X})=\displaystyle\frac{1}{\sqrt{\det(\mathbf{g})}}\partial_i\big(\sqrt{\det(\mathbf{g})}X^i\big).
\]
Substituting $\mathbf{X}=\operatorname{grad}f$ with components $X^i=g^{ij}\partial_j f$, we obtain
\[
\Delta_{\mathbf{g}} f=\operatorname{div}(\operatorname{grad}f)
=\displaystyle\frac{1}{\sqrt{\det(\mathbf{g})}}\partial_i\big(\sqrt{\det(\mathbf{g})}g^{ij}\partial_j f\big).
\]

Thus the local coordinate expression for the Laplace--Beltrami operator is
\[
\Delta_{\mathbf{g}}f=\displaystyle\frac{1}{\sqrt{\det(\mathbf{g})}}
\partial_i\left(\sqrt{\det(\mathbf{g})}g^{ij}\partial_j f\right).
\]
\end{remark}
Combining the Leibniz rules for the gradient and divergence gives the following rule for the Laplacian of the product of two functions:
\begin{proposition}[Leibniz rule for the Laplacian]\label{leibniz para el laplaciano}\index{Leibniz rule for the Laplacian}
 Let $(M,\mathbf{g})$ be a Riemannian manifold with or without boundary. If $f,g\in C^{\infty}(M)$, then \[\Delta_{\mathbf{g}}(fg)=f\Delta_{\mathbf{g}}g+2\langle \operatorname{grad}(f),\operatorname{grad}(g)\rangle_{\mathbf{g}}+g\Delta_{\mathbf{g}}f.\]
\end{proposition}
\begin{proof}
By the two preceding Leibniz rules,
\[
\begin{aligned}
\Delta_{\mathbf g}(fg)
&=\operatorname{div}\bigl(f\operatorname{grad}g
  +g\operatorname{grad}f\bigr)\\
&=f\Delta_{\mathbf g}g
  +\langle\operatorname{grad}f,\operatorname{grad}g\rangle
  +g\Delta_{\mathbf g}f
  +\langle\operatorname{grad}g,\operatorname{grad}f\rangle.
\end{aligned}
\]
Symmetry of \(\mathbf g\) combines the two cross terms.
\end{proof}

\subsection{Integration by parts formulas}

\begin{theorem}[Divergence theorem]\label{teo:divergencia}\index{divergence theorem}
Let $(M,\mathbf{g})$ be a Riemannian manifold with or without boundary, not necessarily orientable,
and let $\mathbf{X}\in \mathfrak{X}(M)$ be a smooth vector field with
compact support. If $\partial M\neq\varnothing$, denote by
$\boldsymbol{\nu}$ the outward unit normal field and by
$\widetilde{\mathbf{g}}$ the induced metric on $\partial M$. Then
\[
\int_{M}\operatorname{div}(\mathbf{X})\,d\lambda_{\mathbf{g}}
= \int_{\partial M}\langle \mathbf{X},\boldsymbol{\nu}\rangle_{\mathbf{g}}
\,d\lambda_{\widetilde{\mathbf{g}}}.
\]
If \(\partial M=\varnothing\), the right-hand side is zero.
\end{theorem}
\begin{proof}
First suppose that \(M\) is oriented and let \(dV_{\mathbf g}\) be the
positive volume form. In an oriented chart, write
\(\mathbf X=X^i\boldsymbol{\partial}_i\) and
\[
dV_{\mathbf g}
=
\sqrt{\det(\mathbf{g})}\,
\mathbf d x^1\wedge\cdots\wedge\mathbf d x^n.
\]
Then
\[
\iota_{\mathbf X}dV_{\mathbf g}
=
\sqrt{\det(\mathbf{g})}
\displaystyle\sum_{i=1}^n(-1)^{i-1}X^i\,
\mathbf d x^1\wedge\cdots\wedge
\widehat{\mathbf d x^i}\wedge\cdots\wedge\mathbf d x^n.
\]
Differentiating each summand, the components of
\(\mathbf d(\sqrt{\det(\mathbf{g})}X^i)\) other than
\(\mathbf d x^i\) vanish when taking the exterior product, because the
corresponding differential already occurs in the \((n-1)\)-form. The
remaining component produces the sign \((-1)^{i-1}\) twice, once from
contraction and once from ordering the differentials. Hence
\[
d\bigl(\iota_{\mathbf X}dV_{\mathbf g}\bigr)
=
\displaystyle\sum_{i=1}^n
\partial_i\!\left(\sqrt{\det(\mathbf{g})}\,X^i\right)
\mathbf d x^1\wedge\cdots\wedge\mathbf d x^n
=
\operatorname{div}(\mathbf X)\,dV_{\mathbf g}.
\]
Apply Stokes' theorem to the compactly supported form
\(\iota_{\mathbf X}dV_{\mathbf g}\). If
\((\mathbf e_1,\dots,\mathbf e_{n-1})\) is a positively oriented orthonormal frame of \(T(\partial M)\), the boundary orientation convention makes
\((\boldsymbol{\nu},\mathbf e_1,\dots,\mathbf e_{n-1})\) positively oriented, and
\[
\begin{aligned}
\bigl(\iota_{\mathbf X}dV_{\mathbf g}\bigr)
(\mathbf e_1,\dots,\mathbf e_{n-1})
&=
dV_{\mathbf g}
(\mathbf X,\mathbf e_1,\dots,\mathbf e_{n-1})\\
&=
\mathbf g(\mathbf X,\boldsymbol{\nu}).
\end{aligned}
\]
Therefore,
\[
\int_M\operatorname{div}(\mathbf X)\,d\lambda_{\mathbf g}
=
\int_{\partial M}
\mathbf g(\mathbf X,\boldsymbol{\nu})\,
d\lambda_{\widetilde{\mathbf g}}
\]
in the orientable case.

Now consider an arbitrary \(M\). The compact set
\(K:=\operatorname{supp}(\mathbf X)\) admits a finite cover by chart domains \(U_1,\dots,U_N\), each of which is orientable. Add \(M\setminus K\) to the cover and choose a subordinate smooth partition of unity. Denote the functions corresponding to \(U_i\) by \(\chi_i\). The function corresponding to \(M\setminus K\) has zero product with \(\mathbf X\), and
\[
\displaystyle\sum_{i=1}^N\chi_i=1,
\qquad
\displaystyle\sum_{i=1}^N\mathbf d\chi_i=0
\]
on a neighborhood of \(K\). Each field \(\chi_i\mathbf X\) has compact
support contained in \(U_i\). Since \(U_i\) is open in a manifold
with boundary, its boundary as a manifold is \(U_i\cap\partial M\);
no artificial boundary appears inside \(M\). The orientable case,
applied to \(U_i\) with either of its two orientations, gives
\[
\int_{U_i}\operatorname{div}(\chi_i\mathbf X)\,d\lambda_{\mathbf g}
=
\int_{\partial M\cap U_i}
\chi_i\mathbf g(\mathbf X,\boldsymbol{\nu})\,
d\lambda_{\widetilde{\mathbf g}}.
\]
The Leibniz rule for divergence implies
\[
\displaystyle\sum_{i=1}^N\operatorname{div}(\chi_i\mathbf X)
=
\left(\displaystyle\sum_{i=1}^N\chi_i\right)\operatorname{div}(\mathbf X)
+
\left\langle
\displaystyle\sum_{i=1}^N\operatorname{grad}\chi_i,\mathbf X
\right\rangle_{\mathbf g}
=
\operatorname{div}(\mathbf X).
\]
Summing the local identities and again using
\(\displaystyle \sum_{i=1}^{N}\chi_i=1\) on \(K\cap\partial M\) gives the global formula. The assertion for empty boundary is the same equality with an integral over the empty set.
\end{proof}

\begin{corollary}[Green identities]
\label{cor:identidades-green-geometricas}
Let $f,h\in C^\infty(M)$ on a Riemannian manifold with or without
boundary, and suppose at least one of them has compact support.
With the convention $\Delta_{\mathbf g}=\operatorname{div}
\operatorname{grad}$ and $\partial_{\boldsymbol\nu}h:=
\mathbf g(\operatorname{grad}h,\boldsymbol\nu)$, we have
\[
\int_M f\Delta_{\mathbf g}h\,d\lambda_{\mathbf g}
=-\int_M\mathbf g(\operatorname{grad}f,\operatorname{grad}h)
\,d\lambda_{\mathbf g}
+\int_{\partial M}f\partial_{\boldsymbol\nu}h
\,d\lambda_{\widetilde{\mathbf g}},
\]
and
\[
\int_M(f\Delta_{\mathbf g}h-h\Delta_{\mathbf g}f)\,d\lambda_{\mathbf g}
=\int_{\partial M}
(f\partial_{\boldsymbol\nu}h-h\partial_{\boldsymbol\nu}f)
\,d\lambda_{\widetilde{\mathbf g}}.
\]
Boundary terms are zero when the boundary is empty.
\end{corollary}
\begin{proof}
The field $f\operatorname{grad}h$ has compact support: if $h$ has
compact support, so does its gradient. The Leibniz rule gives
\[
\operatorname{div}(f\operatorname{grad}h)
=\mathbf g(\operatorname{grad}f,\operatorname{grad}h)
+f\Delta_{\mathbf g}h.
\]
The divergence theorem applied to this field proves the first
identity. Interchanging $f,h$ and subtracting cancels the symmetric
inner product and gives the second.
\end{proof}
\section{Clifford geometry and spin structures}
\label{sec:clifford-spin-dirac}
\index{Clifford algebra}
\index{spin group@$\operatorname{Spin}$ group}
\index{spin structure}

At this point, we have developed Lie groups and Lie algebras, vector bundles,
principal bundles, associated bundles, and principal connections. These are
precisely the geometric ingredients needed to construct
spinor bundles and their connection. We will not yet define an operator
associated with these data: contraction of a Clifford connection will be
studied after developing the general theory of differential
operators; its symbol will be computed only after symbols have been
defined, and its pseudodifferential interpretation is reserved for
the corresponding chapter.

Throughout this section, we use the convention
\[
 vw+wv=-2\mathbf g(v,w)1.
\]
The constructions and conventions may be compared with
\cite[Chapter~I, \S\S1--6; Chapter~II,
\S\S1 and 3--5]{LawsonMichelsohn1989} and
\cite[Chapter~17]{BleeckerBoossIndex}.

\begin{semblanzaHistorica}{Clifford and spinors}
In the nineteenth century, William Kingdon Clifford introduced an algebra combining the
metric with exterior multiplication. Representations of this algebra
give rise to spinors, objects that encode the geometry of the spin
group linearly. On a Riemannian manifold, a spin structure
allows these representations to be globalized and the Levi--Civita
connection to be lifted to the spinor bundle. This geometric construction will underlie
the operators appearing in later chapters.
\end{semblanzaHistorica}

\subsection{Clifford algebras}

\begin{definition}[Clifford algebra]
\label{def:algebra-clifford}
Let $(V,\mathbf g)$ be a real Euclidean space of dimension $n$. Its
\emph{real Clifford algebra} is the unital quotient
\[
 \operatorname{Cl}(V,\mathbf g)
 :=T(V)\big/
 \left\langle v\otimes v+\mathbf g(v,v)1\mid v\in V\right\rangle.
\]
Identify $V$ with its image in the quotient. The complexification is
\[
 \operatorname{Cl}_{\mathbb C}(V,\mathbf g)
 :=\operatorname{Cl}(V,\mathbf g)\otimes_{\mathbb R}\mathbb C.
\]
For the standard Euclidean space, we write
$\operatorname{Cl}_n$ and $\operatorname{Cl}_n^{\mathbb C}$.
\end{definition}

Polarization of $v^2=-|v|_{\mathbf g}^2 1$ gives
\begin{equation}
 vw+wv=-2\mathbf g(v,w)1,
 \qquad v^2=-|v|_{\mathbf g}^2 1.
\label{eq:relacion-clifford}
\end{equation}

\begin{proposition}[Universal property]
\label{prop:propiedad-universal-clifford}
Let $A$ be a real unital associative algebra and let $q\colon V\to A$ be a linear
map such that
\[
 q(v)^2=-|v|_{\mathbf g}^2 1_A
 \qquad(v\in V).
\]
There exists a unique unital algebra homomorphism
\(\widehat q\colon\operatorname{Cl}(V,\mathbf g)\to A\) whose restriction to
$V$ is $q$.
\end{proposition}

\begin{proof}
The universal property of the tensor algebra extends $q$ uniquely
to a unital homomorphism $T(V)\to A$. Every generator
$v\otimes v+|v|^2 1$ of the defining ideal belongs to its kernel. The
homomorphism therefore factors through the quotient. Since
$\operatorname{Cl}(V,\mathbf g)$ is generated by $V$ and the unit, this
factorization is unique.
\end{proof}

\begin{proposition}[Basis, grading, and reversion]
\label{prop:base-graduacion-clifford}
If $(e_1,\ldots,e_n)$ is an orthonormal basis of $V$, then
\[
 1,\qquad e_{i_1}\cdots e_{i_k},
 \qquad 1\leq i_1<\cdots<i_k\leq n,
\]
form a basis of $\operatorname{Cl}(V,\mathbf g)$. In particular,
$\dim_{\mathbb R}\operatorname{Cl}(V,\mathbf g)=2^n$. Moreover:
\begin{enumerate}
\item the automorphism $\alpha$ determined by $\alpha(v)=-v$ induces the
intrinsic grading
\[
 \operatorname{Cl}(V,\mathbf g)
 =\operatorname{Cl}^0(V,\mathbf g)\oplus
 \operatorname{Cl}^1(V,\mathbf g),
\]
given by the eigenspaces of $\alpha$ with eigenvalues $1$ and $-1$;
\item there exists a unique antiautomorphism $\tau$, called \emph{reversion},
such that $\tau(v)=v$ for every $v\in V$, and
\[
 \tau(v_1\cdots v_k)=v_k\cdots v_1.
\]
\end{enumerate}
\end{proposition}

\begin{proof}
The relations $e_i^2=-1$ and $e_ie_j=-e_je_i$ for $i\neq j$ allow us to
order any word and show that the indicated monomials span the
algebra. To prove their independence, define on $\Lambda^\bullet V$
\[
 c_{\Lambda}(v):=\varepsilon(v)-\iota_v,
\]
where $\varepsilon(v)$ is exterior multiplication and $\iota_v$ is
contraction using $\mathbf g$. The identities between contraction and
exterior product give $c_{\Lambda}(v)^2=-|v|^2I$; the universal
property yields a Clifford representation. If $I=(i_1<\cdots<i_k)$,
then
\[
 c_{\Lambda}(e_{i_1})\cdots c_{\Lambda}(e_{i_k})1
 =e_{i_1}\wedge\cdots\wedge e_{i_k}.
\]
The forms on the right constitute a basis of $\Lambda^\bullet V$,
so the Clifford monomials are linearly independent.

The automorphism of $T(V)$ multiplying a tensor of degree $k$ by
$(-1)^k$ preserves the defining ideal and descends to $\alpha$. Its two
eigenspaces give the stated direct sum. Similarly, reversing
the order of factors in $T(V)$ preserves the ideal and descends to the
antiautomorphism $\tau$. In both cases, uniqueness follows because $V$
generates the algebra.
\end{proof}

\subsection{The groups $\operatorname{Pin}(n)$ and $\operatorname{Spin}(n)$}

\begin{definition}[Pin and spin groups]
\label{def:grupos-pin-spin}
Let $S(V)=\{v\in V:|v|_{\mathbf g}=1\}$. Define
\[
 \operatorname{Pin}(V,\mathbf g)
 :=\{v_1\cdots v_k:k\geq0,\ v_j\in S(V)\}
 \subseteq\operatorname{Cl}(V,\mathbf g)^\times
\]
and
\[
 \operatorname{Spin}(V,\mathbf g)
 :=\operatorname{Pin}(V,\mathbf g)\cap
 \operatorname{Cl}^0(V,\mathbf g).
\]
The empty product is $1$. For $V=\mathbb R^n$, write
$\operatorname{Pin}(n)$ and $\operatorname{Spin}(n)$.
\end{definition}

Every unit vector satisfies $v^{-1}=-v$. In particular, the preceding
sets are groups. The parity of a product of unit vectors is
well-defined because the Clifford grading is a direct sum.

\begin{theorem}[Cartan--Dieudonné]
\label{teo:cartan-dieudonne}
Every transformation $A\in\operatorname{O}(V,\mathbf g)$ is a composition of at
most $n$ reflections in hyperplanes through the origin. The
parity of the number of reflections agrees with the sign of $\det A$.
\end{theorem}

\begin{proof}
Proceed by induction on $n$. The result is immediate if $n=0$.
Let $x$ be a unit vector. If $Ax=x$, then $A$ preserves $x^\perp$ and
the induction hypothesis applies to $A|_{x^\perp}$. If $Ax\neq x$, set
\[
 u:=\frac{Ax-x}{|Ax-x|}.
\]
For the reflection $r_u(y)=y-2\mathbf g(y,u)u$, compute
\[
|Ax-x|_{\mathbf g}^2=2(1-\mathbf g(Ax,x)),\qquad
\mathbf g(Ax,Ax-x)=1-\mathbf g(Ax,x)=\tfrac12|Ax-x|_{\mathbf g}^2.
\]
Consequently, $r_u(Ax)=Ax-(Ax-x)=x$. Hence $r_uA$ preserves $x^\perp$ and is a product of at
most $n-1$ reflections in that subspace, extended to $V$ by fixing $x$. Since $r_u^2=I$, this gives a
factorization of $A$ with at most $n$ reflections. Finally, each
reflection has determinant $-1$, which determines the parity.
\end{proof}

For a homogeneous element $a$ of the Clifford algebra, define its twisted adjoint
action on $V$ by
\[
 \widetilde\rho(a)w:=\alpha(a)wa^{-1}.
\]
If $v$ is a unit vector, the Clifford relations give
\[
 \widetilde\rho(v)w=(-v)w(-v)
 =w-2\mathbf g(v,w)v=r_v(w).
\]

\begin{lemma}[Group of units of a finite-dimensional algebra]
\label{lem:unidades-algebra-finita-grupo-lie}
If \(A\) is a finite-dimensional real unital associative algebra, its group
of units \(A^\times\) is open in \(A\) and, with this manifold
structure, is a Lie group.
\end{lemma}

\begin{proof}
For \(a\in A\), let \(L_a\in\operatorname{End}_{\mathbb R}(A)\) be
left multiplication. The element \(a\) is a unit if and only if
\(L_a\) is invertible. Indeed, if \(L_a\) is bijective, there exists \(b\) with
\(ab=1\); moreover,
\(L_a(ba-1)=(ab)a-a=0\), and injectivity gives \(ba=1\). In a basis of
\(A\), the entries of \(L_a\) depend linearly on \(a\), so
\[
 A^\times=\{a\in A:\det L_a\neq0\}
\]
is open. Multiplication is bilinear and, since
\(a^{-1}=L_a^{-1}(1)\), the matrix formula for the inverse proves that
\(a\mapsto a^{-1}\) is smooth.
\end{proof}

\begin{proposition}[The pin group]
\label{prop:compacidad-pin-spin}
For $n\geq1$, the preceding formula defines a surjective homomorphism
\[
 \widetilde\rho\colon\operatorname{Pin}(n)\longrightarrow
 \operatorname{O}(n)
\]
with kernel $\{\pm1\}$. The group $\operatorname{Pin}(n)$ is a compact
subgroup of $\operatorname{Cl}_n^\times$ and hence a Lie group.
$\operatorname{Spin}(n)$ is a compact Lie subgroup.
\end{proposition}

\begin{proof}
Multiplicativity of $\alpha$ shows that $\widetilde\rho$ is a
homomorphism. Its surjectivity follows from
Theorem~\ref{teo:cartan-dieudonne} and the calculation of reflections.

Let $a$ be a homogeneous element of the kernel. If $a$ is even, it commutes with every
vector; if it is odd, it anticommutes with every vector. Write $a$ in the basis
of Proposition~\ref{prop:base-graduacion-clifford}. For a monomial
$e_I$ of degree $k$, we have
\[
 e_Ie_j=(-1)^k e_je_I\quad(j\notin I),
 \qquad
 e_Ie_j=(-1)^{k-1}e_je_I\quad(j\in I).
\]
Comparing coefficients shows that the only even element commuting
with all the $e_j$ is a scalar, and that no nonzero odd element
anticommutes with all of them. Thus, $a=\lambda1$. Since an element
of the kernel is even, it can be written as a product of an even number of
unit vectors; consequently, $a^{-1}=\tau(a)$. Reversion fixes
scalars, so $\lambda^{-1}=\lambda$ and $\lambda=\pm1$.

If $a\in\operatorname{Pin}(n)$, Cartan--Dieudonné gives a product
$b=v_1\cdots v_k$, with $k\leq n$, having the same orthogonal image. The
kernel calculation gives $a=\pm b$. Consequently,
$\operatorname{Pin}(n)$ is the finite union of the images of the compact sets
\[
 \{\pm1\}\times S(\mathbb R^n)^k,
 \qquad 0\leq k\leq n,
\]
under Clifford multiplication. It is compact and, in particular, closed
in the open Lie group $\operatorname{Cl}_n^\times$. The closed subgroup
theorem gives it a Lie group structure. Finally,
$\operatorname{Spin}(n)=\operatorname{Pin}(n)\cap\operatorname{Cl}_n^0$ is
a closed subgroup of $\operatorname{Pin}(n)$ and is also a compact Lie
group.
\end{proof}

On the even part, $\alpha(a)=a$, so the twisted adjoint action is ordinary
conjugation.

\begin{theorem}[Spin covering map]
\label{teo:cubriente-spin}
For $n\geq1$, the formula
\[
 \rho\colon\operatorname{Spin}(n)\longrightarrow\operatorname{SO}(n),
 \qquad \rho(a)v=ava^{-1},
\]
defines a surjective Lie group homomorphism with kernel
$\{\pm1\}$. Its Lie algebra is
\[
 \mathfrak{spin}(n)
 =\operatorname{span}_{\mathbb R}
 \left\{\frac12e_ie_j:1\leq i<j\leq n\right\},
\]
and
\[
 d\rho\left(\frac12e_ie_j\right)e_i=e_j,
 \qquad
 d\rho\left(\frac12e_ie_j\right)e_j=-e_i,
\]
while the action vanishes on the remaining basis vectors. In
particular, $d\rho\colon\mathfrak{spin}(n)\to\mathfrak{so}(n)$ is an
isomorphism and $\rho$ is a two-sheeted covering map:
\[
 1\longrightarrow\{\pm1\}\longrightarrow\operatorname{Spin}(n)
 \xrightarrow{\,\rho\,}\operatorname{SO}(n)\longrightarrow1.
\]
\end{theorem}

\begin{proof}
Proposition~\ref{prop:compacidad-pin-spin} shows that the restriction to
the even part takes values in $\operatorname{SO}(n)$, is surjective, and
has kernel $\{\pm1\}$. For $i<j$, the curve
\[
 \gamma_{ij}(t)=\exp\!\left(\frac t2e_ie_j\right)
 =\cos\frac t2+\sin\frac t2\,e_ie_j
\]
belongs to $\operatorname{Spin}(n)$; indeed,
\[
 \cos\frac t2+\sin\frac t2\,e_ie_j
 =(-e_i)\left(\cos\frac t2\,e_i-\sin\frac t2\,e_j\right)
\]
is a product of two unit vectors. Differentiating
$\rho(\gamma_{ij}(t))v=\gamma_{ij}(t)v\gamma_{ij}(t)^{-1}$ at $t=0$
gives
\[
 d\rho\left(\frac12e_ie_j\right)v
 =\left[\frac12e_ie_j,v\right],
\]
from which the stated formulas follow. Their images form the usual basis
of $\mathfrak{so}(n)$.

On the other hand, $d\rho$ is injective: if $X$ is in its kernel, then
$\rho(\exp(tX))=\exp(t\,d\rho(X))=I$ for every $t$, and discreteness of
$\{\pm1\}$ forces $X=0$. Thus the preceding bivectors span all of
$\mathfrak{spin}(n)$, and $d\rho$ is an isomorphism.

The inverse function theorem gives a neighborhood $W$ of $1$ in
$\operatorname{Spin}(n)$ on which $\rho$ is a diffeomorphism onto a
neighborhood $U$ of $I$. Shrinking $W$ if necessary, $W$ and $-W$ are
disjoint. Since each fiber is a coset of $\{\pm1\}$,
\(\rho^{-1}(U)=W\sqcup(-W)\), and both restrictions are diffeomorphisms
onto $U$. Translations give the same description around
every point of $\operatorname{SO}(n)$. This is precisely the
property of being a two-sheeted covering map.
\end{proof}

\subsection{Complex spinor modules}

The following classification is included to fix the dimensions,
grading, and Hermitian conventions.

\begin{proposition}[Classification of complex Clifford algebras]
\label{prop:clasificacion-clifford-compleja}
For $m\geq0$, there exist complex algebra isomorphisms
\[
 \operatorname{Cl}_{2m}^{\mathbb C}
 \cong\operatorname{End}(\mathbb C^{2^m}),
 \qquad
 \operatorname{Cl}_{2m+1}^{\mathbb C}
 \cong\operatorname{End}(\mathbb C^{2^m})
 \oplus\operatorname{End}(\mathbb C^{2^m}).
\]
\end{proposition}

\begin{proof}
The base cases are
$\operatorname{Cl}_0^{\mathbb C}=\mathbb C$ and
\[
 \operatorname{Cl}_1^{\mathbb C}
 \cong\mathbb C[t]/(t^2+1)
 \cong\mathbb C\oplus\mathbb C.
\]
For the induction step, consider
\[
 J=\begin{pmatrix}1&0\\0&-1\end{pmatrix},\qquad
 A=\begin{pmatrix}0&1\\-1&0\end{pmatrix},\qquad
 B=\begin{pmatrix}0&i\\i&0\end{pmatrix}.
\]
We have $A^2=B^2=-I$, $J^2=I$, and every pair of distinct matrices
anticommutes. If $e_1,\ldots,e_n$ are the generators of
$\operatorname{Cl}_n^{\mathbb C}$, the assignments
\[
 e_j\longmapsto e_j\otimes J,
 \qquad e_{n+1}\longmapsto1\otimes A,
 \qquad e_{n+2}\longmapsto1\otimes B
\]
satisfy the Clifford relations and define a homomorphism
\[
 \operatorname{Cl}_{n+2}^{\mathbb C}\longrightarrow
 \operatorname{Cl}_n^{\mathbb C}\otimes M_2(\mathbb C).
\]
The matrices $I,A,B,J$ generate $M_2(\mathbb C)$; moreover, once
$1\otimes J$ has been obtained, $e_j\otimes1$ can be recovered. Thus the homomorphism is
surjective. Both sides have complex dimension $2^{n+2}$ by
Proposition~\ref{prop:base-graduacion-clifford}, so it is an isomorphism.
Induction gives the two formulas.
\end{proof}

\begin{proposition}[Complex spinor modules]
\label{prop:modulos-espinores-complejos}
If $n=2m$, there exists, up to equivalence, a unique irreducible module
$\Sigma_n$ over $\operatorname{Cl}_n^{\mathbb C}$, and
$\dim_{\mathbb C}\Sigma_n=2^m$. If $n=2m+1$, there are two irreducible
modules of dimension $2^m$; their restrictions to
$\operatorname{Spin}(n)$ are equivalent, and we fix one of the two
Clifford extensions as $\Sigma_n$.

In even dimension, for an oriented orthonormal basis, define
\begin{equation}
 \omega_{\mathbb C}:=i^m c(e_1)\cdots c(e_{2m}).
\label{eq:volumen-complejo-clifford}
\end{equation}
This endomorphism is independent of the chosen basis, satisfies
$\omega_{\mathbb C}^2=I$, and determines the chiral decomposition
\[
 \Sigma_n=\Sigma_n^+\oplus\Sigma_n^-,
 \qquad
 \Sigma_n^\pm=\ker(\omega_{\mathbb C}\mp I).
\]
Multiplication by a vector interchanges the two summands. In every
dimension, a Hermitian inner product can be chosen for which
\begin{equation}
 c(v)^*=-c(v),\qquad v\in\mathbb R^n.
\label{eq:clifford-antiautoadjunta}
\end{equation}
The spin representation is
\[
 \rho_\Sigma\colon\operatorname{Spin}(n)\longrightarrow
 \operatorname{U}(\Sigma_n),\qquad \rho_\Sigma(a)=c(a).
\]
\end{proposition}

\begin{proof}
Proposition~\ref{prop:clasificacion-clifford-compleja} and the
elementary classification of irreducible modules over a matrix
algebra prove the assertions on existence, uniqueness, and dimension. In
odd dimension, the even algebra is identified with the complex Clifford
algebra in dimension $2m$; this is why the two representations restricted
to $\operatorname{Spin}(2m+1)$ are equivalent.

The Clifford relation allows us to reverse the order of the $2m$ factors in
\eqref{eq:volumen-complejo-clifford} and gives
$\omega_{\mathbb C}^2=I$. An oriented orthonormal change of basis multiplies
the alternating product by its determinant, which is one; hence the
element is independent of the basis. Moving a vector past the other
$2m-1$ factors introduces a negative sign, so
$\omega_{\mathbb C}$ anticommutes with $c(v)$, and the latter interchanges the two
eigenspaces.

By Proposition~\ref{prop:compacidad-pin-spin}, we may average
any Hermitian inner product over $\operatorname{Pin}(n)$. For the
resulting inner product, $c(v)$ is unitary when $v$ is a unit vector. Since
$c(v)^{-1}=-c(v)$, we obtain $c(v)^*=-c(v)$; homogeneity extends the
identity to every vector. The restriction to $\operatorname{Spin}(n)$ is
unitary and defines $\rho_\Sigma$.
\end{proof}

Conjugation in the algebra gives the fundamental equivariance identity
\begin{equation}
 \rho_\Sigma(a)c(v)\rho_\Sigma(a)^{-1}
 =c(\rho(a)v),
 \qquad a\in\operatorname{Spin}(n),\ v\in\mathbb R^n.
\label{eq:equivarianza-multiplicacion-clifford-modelo}
\end{equation}

\subsection{The orthonormal frame bundle and spin structures}

Let $(M,\mathbf g)$ be a smooth oriented Riemannian manifold of dimension $n$ with or without
boundary. Define
\[
 \mathbf P_{\operatorname{SO}}(M)_x
 :=\{p\colon\mathbb R^n\to T_xM:
 p\text{ is an orientation-preserving linear isometry}\}.
\]
The right action is $p\cdot A=p\circ A$. A local oriented orthonormal
frame $(\mathbf e_1,\ldots,\mathbf e_n)$ identifies $p$ with the matrix
$A\in\operatorname{SO}(n)$ determined by
$\displaystyle p(e_j^{\mathrm{est}})=\displaystyle\sum_{i=1}^{n} A_i{}^j\mathbf e_i$. These charts make
$\mathbf P_{\operatorname{SO}}(M)\to M$ a principal bundle with group
$\operatorname{SO}(n)$, called the \emph{oriented orthonormal frame
bundle}. Evaluation induces a bundle isomorphism
\begin{equation}
 \mathbf P_{\operatorname{SO}}(M)
 \times_{\operatorname{SO}(n)}\mathbb R^n
 \longrightarrow TM,
 \qquad [p,v]\longmapsto p(v).
\label{eq:haz-tangente-asociado-marcos-so}
\end{equation}
Indeed, the relation $[pA,v]=[p,Av]$ proves that the map is well-defined,
and on each fiber it is a bijective linear isometry.

\begin{definition}[Spinorial structure, or spin structure]
\label{def:estructura-spin}
A \emph{spinorial structure}, also called a \emph{spin structure},
on $(M,\mathbf g)$ consists of a right principal bundle
$\mathbf P_{\operatorname{Spin}}(M)\to M$ with group
$\operatorname{Spin}(n)$ and a smooth map
\[
 \Lambda\colon\mathbf P_{\operatorname{Spin}}(M)\longrightarrow
 \mathbf P_{\operatorname{SO}}(M)
\]
covering $\operatorname{id}_M$ and satisfying
\begin{equation}
 \Lambda(p\cdot a)=\Lambda(p)\cdot\rho(a).
\label{eq:equivarianza-estructura-spin}
\end{equation}
Two structures are equivalent if there exists a principal bundle
isomorphism commuting with their maps $\Lambda$. An oriented Riemannian
manifold is called \emph{spin} if it admits a spin structure; when
speaking of its spinors, an equivalence class of such structures
is assumed to have been chosen.
\end{definition}

In principal trivializations, the identity
\eqref{eq:equivarianza-estructura-spin} shows that $\Lambda$ has the form
\((x,a)\mapsto(x,b(x)\rho(a))\). By
Theorem~\ref{teo:cubriente-spin}, $\Lambda$ is a two-sheeted covering
map restricting on each fiber to the spin map. Existence of these data is
not automatic. Appendix~\ref{ap:topologia}, after introducing the necessary
characteristic classes, proves the topological existence
criterion and classifies spin structures.

\begin{definition}[Spinor bundle and Clifford multiplication]
\label{def:haz-espinores}
Let $(M,\mathbf g)$ be a spin manifold with a chosen structure. Its
\emph{complex spinor bundle} is the associated bundle
\[
 \boldsymbol\Sigma M
 :=\mathbf P_{\operatorname{Spin}}(M)
 \times_{\rho_\Sigma}\Sigma_n.
\]
For $\xi\in T_x^*M$, $p\in\mathbf P_{\operatorname{Spin}}(M)_x$, and
$\psi\in\Sigma_n$, let
\(v=\Lambda(p)^{-1}(\xi^\sharp)\in\mathbb R^n\).
\emph{Clifford multiplication} is defined by
\begin{equation}
 c(\xi)[p,\psi]:=[p,c(v)\psi].
\label{eq:definicion-multiplicacion-clifford-espinores}
\end{equation}
\end{definition}

The definition is independent of the representative. Indeed, replacing
$p$ by $pa$ changes the coordinates of $\xi^\sharp$ and the spinor
to $\rho(a)^{-1}v$ and $\rho_\Sigma(a)^{-1}\psi$; the identity
\eqref{eq:equivarianza-multiplicacion-clifford-modelo} gives the same
class. Therefore, \eqref{eq:definicion-multiplicacion-clifford-espinores}
defines a smooth bundle homomorphism
\[
 c\colon\operatorname{Cl}_{\mathbb C}(T^*M,\mathbf g^{-1})
 \longrightarrow\operatorname{End}(\boldsymbol\Sigma M)
\]
and, for real covectors,
\begin{equation}
 c(\xi)c(\eta)+c(\eta)c(\xi)
 =-2\mathbf g^{-1}(\xi,\eta)I,
 \qquad c(\xi)^*=-c(\xi).
\label{eq:clifford-haz-espinores}
\end{equation}
If $n=2m$, chirality is invariant under $\operatorname{Spin}(n)$ and
gives
\[
 \boldsymbol\Sigma M=\boldsymbol\Sigma^+M\oplus
 \boldsymbol\Sigma^-M;
\]
multiplication by a covector is odd with respect to this grading.

\subsection{The spin connection}

Denote by $\boldsymbol\omega^{\operatorname{SO}}$ the Levi--Civita principal connection
form on
$\mathbf P_{\operatorname{SO}}(M)$. Since $d\rho$ is an isomorphism,
we may define
\begin{equation}
 \boldsymbol\omega^{\operatorname{Spin}}
 :=(d\rho)^{-1}\Lambda^*
 \boldsymbol\omega^{\operatorname{SO}}.
\label{eq:levantamiento-conexion-spin-principal}
\end{equation}

\begin{proposition}[Spin connection]
\label{prop:compatibilidad-conexion-spinorial}
The form \eqref{eq:levantamiento-conexion-spin-principal} is the unique
principal connection on $\mathbf P_{\operatorname{Spin}}(M)$ that projects
to the Levi--Civita connection. The associated connection
$\nabla^{\boldsymbol\Sigma}$ is Hermitian, preserves the grading in
even dimension, and satisfies
\begin{equation}
 \nabla_{\mathbf X}^{\boldsymbol\Sigma}
 \bigl(c(\xi)\psi\bigr)
 =c(\nabla_{\mathbf X}\xi)\psi
 +c(\xi)\nabla_{\mathbf X}^{\boldsymbol\Sigma}\psi.
\label{eq:compatibilidad-conexion-clifford-spin}
\end{equation}
If $(\mathbf e_1,\ldots,\mathbf e_n)$ is an oriented orthonormal frame admitting
a local spin lift and
\[
 \omega_{ij}(\mathbf X)
 :=\mathbf g(\nabla_{\mathbf X}\mathbf e_i,\mathbf e_j),
\]
then, in this trivialization,
\begin{equation}
 \bigl(\nabla_{\mathbf X}^{\boldsymbol\Sigma}\psi\bigr)^\wedge
 =\mathbf X(\widehat\psi)
 +\frac14\displaystyle\sum_{i,j=1}^n
 \omega_{ij}(\mathbf X)c(e_i)c(e_j)\widehat\psi.
\label{eq:formula-local-conexion-spinorial}
\end{equation}
\end{proposition}

\begin{proof}
The pullback $\Lambda^*\boldsymbol\omega^{\operatorname{SO}}$ reproduces
vertical fields after applying $d\rho$ and is equivariant under the
adjoint action. Since $d\rho$ intertwines the adjoint representations, its
composition with $(d\rho)^{-1}$ satisfies the two axioms of a principal connection
form. Any connection projecting to Levi--Civita must
have this form, which proves uniqueness.

The calculation in Theorem~\ref{teo:cubriente-spin} identifies the element
of $\mathfrak{spin}(n)$ corresponding to the skew-symmetric matrix
$(\omega_{ij})$ with
\(\displaystyle \frac14\displaystyle\sum_{i,j=1}^{n}\omega_{ij}e_ie_j\). The local formula for associated
connections gives \eqref{eq:formula-local-conexion-spinorial}. Each
$c(e_i)c(e_j)$, with $i\neq j$, is skew-adjoint, so the connection
is Hermitian. The complex volume element is invariant under
$\operatorname{Spin}(n)$ and parallel; hence the grading is preserved.
Finally, differentiating the equivariance identity
\eqref{eq:equivarianza-multiplicacion-clifford-modelo}, or substituting
the local formula directly, gives
\eqref{eq:compatibilidad-conexion-clifford-spin}.
\end{proof}

\subsection{Clifford modules and connections}

The later analytic construction requires only the structure
just obtained, without necessarily requiring a global spin structure.

\begin{definition}[Hermitian Clifford module]
\label{def:modulo-conexion-clifford}
A \emph{Hermitian Clifford module} over $(M,\mathbf g)$ is a complex
Hermitian bundle $\mathbf S\to M$ with a smooth unital homomorphism of
algebra bundles
\[
 c\colon\operatorname{Cl}_{\mathbb C}(T^*M,\mathbf g^{-1})
 \longrightarrow\operatorname{End}(\mathbf S)
\]
such that $c(\xi)^*=-c(\xi)$ for every real covector $\xi$. Equivalently,
it suffices to give a smooth linear map $c\colon T^*M\to\operatorname{End}(\mathbf S)$
satisfying \eqref{eq:clifford-haz-espinores}; the equivalence is the
Clifford universal property applied fiberwise.

A Hermitian connection $\nabla^{\mathbf S}$ is a \emph{Clifford
connection} if
\begin{equation}
 \nabla_{\mathbf X}^{\mathbf S}(c(\xi)\mathbf s)
 =c(\nabla_{\mathbf X}\xi)\mathbf s
 +c(\xi)\nabla_{\mathbf X}^{\mathbf S}\mathbf s.
\label{eq:definicion-conexion-clifford}
\end{equation}
If $\mathbf S=\mathbf S^+\oplus\mathbf S^-$ is graded, we also require
the connection to preserve the grading and $c(\xi)$ to be odd.
\end{definition}

The bundle $\boldsymbol\Sigma M$ with
$\nabla^{\boldsymbol\Sigma}$ is the fundamental spinorial example.

\begin{example}[Differential forms]
\label{ej:modulo-clifford-formas}
On $\Lambda^kT^*M$ we use the determinant inner product: for covectors,
\[
\left\langle\alpha_1\wedge\cdots\wedge\alpha_k,
\beta_1\wedge\cdots\wedge\beta_k\right\rangle_{\Lambda^k}
=\det\bigl(\langle\alpha_i,\beta_j\rangle_{\mathbf g^{-1}}\bigr)_{i,j=1}^k.
\]
Equivalently, this is $1/k!$ times the inner product obtained by viewing
alternating forms as covariant tensors. Thus increasing exterior
monomials in an orthonormal coframe are orthonormal. We use the direct
sum of these products in different degrees and their Hermitian
extensions, linear in the first argument, on complex forms. This is the
convention for all exterior bundles, their $L^2(M,\Lambda^kT^*M)$ products, and the Hodge
operators below.

Every Riemannian manifold, whether spin or not, has the graded module
\[
 \mathbf S=\Lambda^\bullet T^*M\otimes\mathbb C,
 \qquad c(\xi)=\varepsilon(\xi)-\iota_{\xi^\sharp}.
\]
The identities
\[
 \iota_{\xi^\sharp}\varepsilon(\eta)
 +\varepsilon(\eta)\iota_{\xi^\sharp}
 =\mathbf g^{-1}(\xi,\eta)I
\]
prove the Clifford relation. For the Hermitian inner product induced by
$\mathbf g$, we have
$\varepsilon(\xi)^*=\iota_{\xi^\sharp}$, so $c(\xi)^*=-c(\xi)$. The
connection induced by Levi--Civita commutes with exterior product and
contraction in the sense of the Leibniz rule; it is therefore a
Clifford connection. The grading is the parity of forms.
\end{example}

\subsection{$\operatorname{Spin}^c$ structures}

\begin{definition}[The group $\operatorname{Spin}^c(n)$]
\label{def:grupo-spinc}
Define the compact Lie group
\[
 \operatorname{Spin}^c(n)
 :=\bigl(\operatorname{Spin}(n)\times\operatorname{U}(1)\bigr)
 /\{(1,1),(-1,-1)\}.
\]
The homomorphisms
\[
 \rho_c[a,z]=\rho(a),
 \qquad \det_c[a,z]=z^2
\]
take values in $\operatorname{SO}(n)$ and $\operatorname{U}(1)$,
respectively. The complex spin representation is
\begin{equation}
 \rho_\Sigma^c[a,z]=z\rho_\Sigma(a).
\label{eq:representacion-spinorial-spinc}
\end{equation}
\end{definition}

All three formulas are well-defined on the quotient: replacing $(a,z)$
by $(-a,-z)$ changes none of them, since
$\rho(-a)=\rho(a)$ and $\rho_\Sigma(-1)=-I$. The quotient by the finite central
subgroup inherits a Lie group structure. Indeed, the action
of \(\{(1,1),(-1,-1)\}\) by translations is free and, since the group is
finite, properly discontinuous. A sufficiently small neighborhood of
each point is disjoint from its nontrivial translate; its projection to the quotient
is a chart, and in these charts the induced multiplication and inversion
are smooth.

\begin{definition}[$\operatorname{Spin}^c$ structure]
\label{def:estructura-spinc-geometrica}
A \emph{$\operatorname{Spin}^c$ structure} on $(M,\mathbf g)$ consists of
a principal bundle
$\mathbf P_{\operatorname{Spin}^c}(M)\to M$ and a map
\[
 \Lambda_c\colon\mathbf P_{\operatorname{Spin}^c}(M)
 \longrightarrow\mathbf P_{\operatorname{SO}}(M)
\]
covering the identity and satisfying
$\Lambda_c(pa)=\Lambda_c(p)\rho_c(a)$. Its \emph{determinant line} and
\emph{spinor bundle} are, respectively,
\[
 \mathbf L
 :=\mathbf P_{\operatorname{Spin}^c}(M)\times_{\det_c}\mathbb C,
 \qquad
 \boldsymbol\Sigma^cM
 :=\mathbf P_{\operatorname{Spin}^c}(M)
 \times_{\rho_\Sigma^c}\Sigma_n.
\]
\end{definition}

The equivariance identity
\eqref{eq:equivarianza-multiplicacion-clifford-modelo} proves, exactly
as in the spin case, that
\[
 c(\xi)[p,\psi]
 :=[p,c(\Lambda_c(p)^{-1}\xi^\sharp)\psi]
\]
defines Clifford multiplication on $\boldsymbol\Sigma^cM$. In
even dimension, it gives the grading
$\boldsymbol\Sigma^{c,+}M\oplus\boldsymbol\Sigma^{c,-}M$.

A unitary connection on $\mathbf L$, together with the Levi--Civita
connection, uniquely determines a principal connection on
$\mathbf P_{\operatorname{Spin}^c}(M)$. Indeed, the Lie algebra map
\[
 (d\rho_c,d\det_c)\colon
 \mathfrak{spin}^c(n)\longrightarrow
 \mathfrak{so}(n)\oplus i\mathbb R
\]
is an isomorphism: on the first factor it is $d\rho$, and on the second,
the differential of $z\mapsto z^2$ is multiplication by $2$. The pullbacks of
the two connection forms, followed by the inverse of this isomorphism,
give the desired principal connection. The associated connection on
$\boldsymbol\Sigma^cM$ is Hermitian and Clifford by the same argument
as in Proposition~\ref{prop:compatibilidad-conexion-spinorial}.

The topological existence criterion is stated and proved, after the
required characteristic classes have been introduced, in
Theorem~\ref{teo:spinc-apendice} of Appendix~\ref{ap:topologia}.

We have thus constructed the purely geometric data: the spin group and its
covering map, spin structures, spinor bundles, Clifford
multiplication, and their connections. In
Chapter~\ref{cap:operadores-diferenciales-parciales-haces}, a Clifford
connection will be contracted to obtain a differential operator. Its symbol will be
studied in Chapter~\ref{cap:simbolo-principal-operadores-elipticos}, and
its pseudodifferential symbol, with the appropriate Fourier convention, only
in Chapter~\ref{cap:operadores-pseudodiferenciales}.

\part{Nonlinear analysis in Euclidean spaces}
\chapter*{Introduction to Part II}
\addcontentsline{toc}{chapter}{Introduction to Part II}
\markboth{Part II. Nonlinear analysis in Euclidean spaces}{Introduction to Part II}

When an estimate on a manifold is proved within a chart, we are ultimately using a result formulated on an open subset of $\mathbb R^n$. Before returning to geometry, it is therefore useful to develop the Euclidean model carefully and determine which tools depend on the linear structure of the space, which survive a change of coordinates, and which can later be localized using cutoff functions.

Sobolev spaces arise when derivatives are allowed to be understood in the weak sense. Distributions extend this idea further, making it possible to treat singularities and data that are not ordinary functions. The Fourier transform allows differentiation and regularity to be studied from the frequency viewpoint. Interpolation theory then provides a common language for comparing different levels of regularity and leads naturally to the Bessel, Slobodeckij, Besov, and Triebel--Lizorkin scales. Orders of regularity are thus no longer restricted to integers.

This part supplies many of the tools that will later be used with almost no change to their mathematical substance: Sobolev inequalities, Rellich--Kondrashov compactness, multipliers, changes of variables, extension operators, convolutions, and localization arguments. On returning to a manifold, what changes is how these tools are transported and how the constants are controlled. The Euclidean work in these chapters makes it possible to distinguish clearly what belongs to analysis from what geometry must provide.

\begin{semblanzaHistorica}{From classical solutions to weak solutions}
For a long time, solving a differential equation meant finding a sufficiently differentiable function that satisfied it pointwise. Twentieth-century analysis decisively broadened this idea. Sobolev incorporated weak derivatives into normed spaces; Schwartz constructed a stable theory of distributions; Calderón and Zygmund developed fundamental tools of harmonic analysis; Gagliardo and Slobodeckij studied noninteger regularity; and Lions, Peetre, and Triebel organized broad scales of function spaces. These theories did not replace classical solutions: they made it possible to begin with less regularity and understand the conditions under which it could later be recovered.
\end{semblanzaHistorica}

\chapter{Sobolev spaces on open subsets of \texorpdfstring{$\mathbb{R}^{n}$}{Rn}}
\label{cap:sobolev-euclidiano}

Sobolev spaces arise naturally in the Dirichlet problem for the Poisson equation. Let $\Omega\subseteq\mathbb R^n$ be an open subset and consider
\[
-\Delta u=f\quad\text{on }\Omega,
\qquad
u=0\quad\text{on }\partial\Omega.
\]
Suppose for the moment that $u$ and $\varphi$ are smooth, that $u$ satisfies the boundary condition, and that $\varphi$ is a test function vanishing on the boundary. Integration by parts gives
\[
\int_\Omega \nabla u\cdot\nabla\varphi\,dx
=
\int_\Omega f\varphi\,dx.
\]
This identity does not require classical second derivatives of $u$: it is enough for $u$ and its gradient to have suitable integrability. This leads to weak derivatives, to $W^{1,2}(\Omega)$, and to the subspace $W_0^{1,2}(\Omega)$ that encodes the Dirichlet condition. The trace provides an interpretation of boundary values, while the Poincaré inequality, embeddings, and compactness theorems provide the estimates needed to solve the problem.

The same equation is obtained as the critical point condition for the functional
\[
\mathcal E(u)
=
\frac{1}{2}\int_\Omega |\nabla u|^2\,dx
-\int_\Omega fu\,dx.
\]
Indeed, its Euler--Lagrange equation is precisely the preceding weak identity. This relationship between a PDE and an energy will serve as the model for the variational methods developed at the end of the book.

The nonlinear problem
\[
-\operatorname{div}\bigl(|\nabla u|^{p-2}\nabla u\bigr)=f
\]
is formally associated with
\[
\mathcal E_p(u)
=
\frac{1}{p}\int_\Omega |\nabla u|^p\,dx
-\int_\Omega fu\,dx.
\]
This example leads to $W^{1,p}(\Omega)$ and shows why convexity and reflexivity matter in the study of minimizing sequences. The theory of the $p$-Laplacian will not be developed here; it will guide the definitions and estimates in this chapter.

\begin{semblanzaHistorica}{Sobolev and the shift toward weak regularity}
In the 1930s, Sergei Sobolev showed that the regularity useful for an equation need not coincide with classical pointwise differentiability. By incorporating weak derivatives into normed spaces, he made it possible to combine integration by parts, completeness, and regularity estimates within a single framework. His embedding theorems also explain a central phenomenon in analysis: a certain degree of integrability of the derivatives forces a function to gain integrability, continuity, or even Hölder regularity. This chapter develops these ideas in the Euclidean setting before transferring them to manifolds and vector bundles~\cite{Sobolev1938}.
\end{semblanzaHistorica}

\section{Weak derivatives}

We begin with the elementary formula that motivates the definition.

\begin{proposition}[Integration by parts]\label{porpartes}\index{integration by parts}
If $f\in C^1(\Omega)$ and $\phi\in C_c^1(\Omega)$, then
\[
\int_\Omega \frac{\partial f}{\partial x_i}\phi
+
\int_\Omega f\frac{\partial\phi}{\partial x_i}
=0
\]
for every $i\in\{1,\dots,n\}$.
\end{proposition}
\begin{proof}
The extension of $\phi$ by zero belongs to $C_c^1(\mathbb R^n)$, since the support of $\phi$ is a compact subset of $\Omega$. We may therefore work in $\mathbb R^n$. Choose $M>0$ such that $\supp(\phi)\subset(-M,M)^n$. For $i=1$ and each $\widehat x=(x_2,\dots,x_n)\in\mathbb R^{n-1}$, Theorem~\ref{teo:b5-fundamental-calculo-riemann-banach} gives
\[
\int_{-M}^{M}\frac{\partial\phi}{\partial x_1}(x_1,\widehat x)\,dx_1
=
\phi(M,\widehat x)-\phi(-M,\widehat x)
=0.
\]
Theorem~\ref{Fubini} allows us to integrate with respect to the remaining variables and conclude that $\displaystyle\int_\Omega \displaystyle\frac{\partial\phi}{\partial x_1}=0$. The other indices are treated in the same way. Applying this identity to $f\phi\in C_c^1(\Omega)$ and using the product rule, we obtain
\[
0
=
\int_\Omega\frac{\partial(f\phi)}{\partial x_i}
=
\int_\Omega\frac{\partial f}{\partial x_i}\phi
+
\int_\Omega f\frac{\partial\phi}{\partial x_i}.
\]
\end{proof}
The preceding identity suggests the following definition.
\begin{definition}\label{def:sobolev-en-abiertos-euclidianos-debilmente-diferenciable-en}\index{weakly differentiable}
Let $u\in L_{\text{loc}}^{1}(\Omega)$. We say that $u$ is \textbf{weakly differentiable on $\Omega$} if there exist $v_{1},\dots,v_{n}\in L_{\text{loc}}^{1}(\Omega)$ such that

\[\int_{\Omega}u\frac{\partial\varphi}{\partial x_{i}}+\int_{\Omega}v_{i}\varphi=0\] for every $\varphi\in C_{c}^{\infty}(\Omega)$ and every $i\in \{1,\dots,n\}$. Each $v_i$ is called the \textbf{$i$th weak derivative of $u$ on $\Omega$} and is denoted by $D_i u$.
\end{definition}

Differential operators constructed from partial derivatives thus admit a weak version. In particular, one can define the weak gradient, the weak Laplacian, and, analogously, higher-order operators.

\begin{remark}\label{obs:sobolev-en-abiertos-euclidianos-satisfacen-c-t-p-proposicion-prop}
If $v,w\in L^1_{\mathrm{loc}}(\Omega)$ satisfy
\[
\int_{\Omega} u\frac{\partial\phi}{\partial x_i} = -\int_{\Omega} v\phi
\quad\text{and}\quad
\int_{\Omega} u\frac{\partial\phi}{\partial x_i} = -\int_{\Omega} w\phi
\qquad \forall\phi\in C_c^\infty(\Omega),
\]
then $v=w$ a.e.\ (by Proposition~\ref{prop unicidad derivada debil prop 14.49}).
\end{remark}

The weak derivative agrees with the classical derivative whenever the latter exists. If $u\in C^1(\Omega)$, Proposition~\ref{porpartes} gives \[\int_{\Omega}\frac{\partial u}{\partial x_{i}}\phi+\int_{\Omega}u\frac{\partial \phi}{\partial x_{i}}=0\] so that $u$ is weakly differentiable and \[D_{i}u=\frac{\partial u}{\partial x_{i}} \text{ $\forall i\in\{1,\dots,n\}$}.\]

Weak derivatives have the following elementary properties.
\begin{proposition}\label{linealidad derivada debil}
 \begin{enumerate}[label=(\alph*)]
 \item If $u,v\in L_{\text{loc}}^{1}(\Omega)$ are weakly differentiable and $\lambda,\mu\in \mathbb{R}$, then $\lambda u+\mu v$ is weakly differentiable and \[D_{i}(\lambda u+\mu v)=\lambda D_{i}u+\mu D_{i}v\] for each $i\in \{1,\dots,n\}$.
 \item If $u\in L_{\text{loc}}^{1}(\Omega)$ is weakly differentiable on $\Omega$ and $\zeta\in C_{c}^{\infty}(\Omega)$, then the extension of $\zeta u$ by zero is weakly differentiable on $\mathbb{R}^{n}$ and \[D_{i}(\zeta u)=\frac{\partial \zeta}{\partial x_{i}}u+\zeta D_{i}u\]
 \end{enumerate}
\end{proposition}
\begin{proof}
 \begin{enumerate}[label=(\alph*)]
 \item For each $i\in \{1,\dots,n\}$ and $\phi\in C_{c}^{\infty}(\Omega)$, \[\int_{\Omega}\left[ (\lambda u+\mu v)\frac{\partial \phi}{\partial x_{i}}+(\lambda D_{i}u+\mu D_{i}v)\phi\right]=\lambda\int_{\Omega}\left[u\frac{\partial \phi}{\partial x_{i}}+(D_{i}u)\phi\right]+\mu\int_{\Omega}\left[v\frac{\partial \phi}{\partial x_{i}}+(D_{i}v)\phi\right]=0\] where the second equality uses the weak differentiability of $u$ and $v$. This proves that $\lambda u+\mu v$ is weakly differentiable and $D_{i}(\lambda u+\mu v)=\lambda D_{i}u+\mu D_{i}v$.
 \item Let $\phi\in C_{c}^{\infty}(\mathbb{R}^{n})$. Then $\zeta\phi\in C_{c}^{\infty}(\Omega)$ and \[\int_{\mathbb{R}^{n}}\zeta u\frac{\partial\phi}{\partial x_{i}}+\int_{\mathbb{R}^{n}}\left(\zeta D_{i}u+u\frac{\partial\zeta}{\partial x_{i}}\right)\phi=\int_{\Omega}u\frac{\partial(\zeta\phi)}{\partial x_{i}}+\int_{\Omega}(D_{i}u)(\zeta\phi)=0.\]
 \end{enumerate}
\end{proof}
A function whose weak derivatives vanish retains the constancy property of differentiable functions with zero differential.

\begin{lemma}[Functions with vanishing weak derivatives]
\label{lem:derivadas-debiles-nulas-euclidiano}
\index{function with vanishing weak derivatives}
Let $\Omega\subseteq\mathbb R^n$ be an open subset and let $v\in L^1_{\mathrm{loc}}(\Omega)$ be weakly differentiable. If
\begin{equation}
\label{eq:derivadas-debiles-euclidianas-nulas-constancia-cap9}
D_jv=0
\qquad\text{almost everywhere on }\Omega,\quad
j\in\{1,\dots,n\},
\end{equation}
then $v$ is constant almost everywhere on each connected component of $\Omega$.
\end{lemma}

\begin{proof}
It suffices to consider the case in which $\Omega$ is connected and nonempty. Fix a ball $B=B_{\mathrm{euc}}(x_0,r)$ with $\overline B\subseteq\Omega$. Choose $R>r$ such that $\overline{B_{\mathrm{euc}}(x_0,R)}\subseteq\Omega$. Let $\rho$ be the mollifier from Theorem~\ref{teo:molificadores-euclidianos} and, for $0<\varepsilon<R-r$, define
\[
v_\varepsilon(x)
:=\int_\Omega v(y)\rho_\varepsilon(x-y)\,dy,
\qquad x\in B.
\]
For each $x\in B$, the function $y\mapsto\rho_\varepsilon(x-y)$ has compact support in $B_{\mathrm{euc}}(x_0,R)$. The local integrability of $v$ and part~(a) of that theorem allow differentiation under the integral. Applying the definition of weak derivative to this test function, we obtain
\[
\begin{aligned}
\partial_jv_\varepsilon(x)
&=\int_\Omega v(y)\partial_{x_j}\rho_\varepsilon(x-y)\,dy\\
&=-\int_\Omega v(y)\partial_{y_j}\rho_\varepsilon(x-y)\,dy\\
&=\int_\Omega D_jv(y)\rho_\varepsilon(x-y)\,dy
=0.
\end{aligned}
\]
Therefore, $v_\varepsilon$ is smooth and has zero differential on the connected ball $B$, so that $v_\varepsilon=c_\varepsilon$ on $B$ for some constant $c_\varepsilon\in\mathbb R$.

Part~(b) of Theorem~\ref{teo:molificadores-euclidianos}, applied with $p=1$ and the compact set $\overline B$, gives $v_\varepsilon\to v$ in $L^1(B)$. Set
\[
c_B:=\frac{1}{\lambda_n(B)}\int_B v\,d\lambda_n.
\]
Since $v_\varepsilon=c_\varepsilon$ on $B$,
\[
|c_\varepsilon-c_B|
\leq\frac{1}{\lambda_n(B)}
\|v_\varepsilon-v\|_{L^1(B)}\longrightarrow0.
\]
Consequently,
\[
\|v-c_B\|_{L^1(B)}
\leq\|v-v_\varepsilon\|_{L^1(B)}
+\lambda_n(B)|c_\varepsilon-c_B|
\longrightarrow0,
\]
and $v=c_B$ almost everywhere on $B$.

If two relatively compact balls in $\Omega$ intersect, their intersection is a nonempty open set and has positive measure. The corresponding constants agree, since both equal $v$ almost everywhere on that intersection. Thus, we can define $c(x):=c_B$ using any ball $B\Subset\Omega$ containing $x$. This function is well defined and locally constant. Since $\Omega$ is connected, it has a single value $c\in\mathbb R$.

Finally, second countability of $\Omega$ allows us to choose a countable cover by such balls. On each ball, $v=c$ outside a set of measure zero; the countable union of these sets also has measure zero. Therefore, $v=c$ almost everywhere on $\Omega$.
\end{proof}

\section{Sobolev spaces}

The notion of weak derivative allows us to define Sobolev spaces on open subsets of $\mathbb R^n$.
\begin{definition}\label{def:sobolev-en-abiertos-euclidianos-espacio-de-sobolev}\index{Sobolev space}\glsadd{espacio-sobolev}
Let $p\in [1,\infty]$. The \textbf{Sobolev space} $W^{1,p}(\Omega)$ is defined by \[W^{1,p}(\Omega):=\{u\in L^{p}(\Omega)\mid u \text{ is weakly differentiable on $\Omega$ and } D_{i}u\in L^{p}(\Omega)\text{ for every $i\in \{1,\dots,n\}$}\}.\]
\end{definition}
We introduce the normed structure as follows.
\begin{definition}\label{def:sobolev-en-abiertos-euclidianos-norma}\index{norm of the first-order Sobolev space}
Let $p\in[1,\infty]$. For $u\in W^{1,p}(\Omega)$, define
\[\|u\|_{W^{1,p}(\Omega)}:=\begin{cases}
 (\|u\|_{L^{p}(\Omega)}^{p}+\displaystyle\sum_{i=1}^{n}\|D_{i}u\|_{L^{p}(\Omega)}^{p})^{\frac{1}{p}} & \text{if $p\in[1,\infty)$}\\
 \max\{\|u\|_{L^{\infty}(\Omega)},\|D_{1}u\|_{L^{\infty}(\Omega)},\dots,\|D_{n}u\|_{L^{\infty}(\Omega)}\} & \text{if $p=\infty$}
 \end{cases}\]
\end{definition}

\subsection{Examples: corners, jumps, and singularities}

Sobolev regularity requires an integration-by-parts identity and an
integrability condition. It does not require classical derivatives at
every point. Nor does computing a classical derivative away from a set
of measure zero suffice: the weak identity may detect a jump across that set.

\begin{example}[Powers with a corner on an interval]
\label{ej:sobolev-potencias-intervalo}
Let $I:=(-1,1)$, $\beta>0$, and $u_\beta(t):=|t|^\beta$.
The function is continuous and smooth away from $0$. Its weak derivative is
\[
 v_\beta(t):=\beta\,\operatorname{sgn}(t)|t|^{\beta-1}
 \quad(t\neq0),\qquad v_\beta(0):=0.
\]
Indeed, $v_\beta$ is locally integrable because $\beta-1>-1$.
For $\varphi\in C_c^\infty(I)$ and $0<\varepsilon<1$, integration by
parts on the two separated intervals gives
\[
 \int_{-1}^{-\varepsilon}u_\beta\varphi'
 +\int_\varepsilon^1u_\beta\varphi'
 =-\int_{(-1,-\varepsilon)\cup(\varepsilon,1)}v_\beta\varphi
   +\varepsilon^\beta\bigl(\varphi(-\varepsilon)-\varphi(\varepsilon)\bigr).
\]
The boundary term satisfies
\[
 \left|\varepsilon^\beta
       \bigl(\varphi(-\varepsilon)-\varphi(\varepsilon)\bigr)\right|
 \leq 2\varepsilon^\beta\|\varphi\|_{L^\infty(I)}.
\]
The integrals over $(-\varepsilon,\varepsilon)$ satisfy the estimates
\begin{align*}
 \left|\int_{-\varepsilon}^{\varepsilon}
          u_\beta(t)\varphi'(t)\,dt\right|
 &\leq \|\varphi'\|_{L^\infty(I)}
        \int_{-\varepsilon}^{\varepsilon}|t|^\beta\,dt
  =\frac{2\varepsilon^{\beta+1}}{\beta+1}
       \|\varphi'\|_{L^\infty(I)},\\
 \left|\int_{-\varepsilon}^{\varepsilon}
          v_\beta(t)\varphi(t)\,dt\right|
 &\leq \beta\|\varphi\|_{L^\infty(I)}
        \int_{-\varepsilon}^{\varepsilon}|t|^{\beta-1}\,dt
  =2\varepsilon^\beta\|\varphi\|_{L^\infty(I)}.
\end{align*}
Adding the integrals over the three intervals gives
\begin{align*}
 \int_Iu_\beta(t)\varphi'(t)\,dt
 +\int_Iv_\beta(t)\varphi(t)\,dt
 &=\varepsilon^\beta
      \bigl(\varphi(-\varepsilon)-\varphi(\varepsilon)\bigr)\\
 &\quad+\int_{-\varepsilon}^{\varepsilon}
            u_\beta(t)\varphi'(t)\,dt
       +\int_{-\varepsilon}^{\varepsilon}
            v_\beta(t)\varphi(t)\,dt.
\end{align*}
The preceding estimates imply
\[
 \left|\int_Iu_\beta(t)\varphi'(t)\,dt
       +\int_Iv_\beta(t)\varphi(t)\,dt\right|
 \leq 4\varepsilon^\beta\|\varphi\|_{L^\infty(I)}
      +\frac{2\varepsilon^{\beta+1}}{\beta+1}
         \|\varphi'\|_{L^\infty(I)}.
\]
Since $\beta>0$, the right-hand side tends to zero as
$\varepsilon\to0^+$. The left-hand side is independent of $\varepsilon$;
therefore,
\[
 \int_Iu_\beta(t)\varphi'(t)\,dt
 =-\int_Iv_\beta(t)\varphi(t)\,dt.
\]

For $1\leq p<\infty$,
\[
 \|u_\beta\|_{L^p(I)}^p=\frac{2}{\beta p+1},
 \qquad
 \int_I|v_\beta|^p
 =2\beta^p\int_0^1t^{p(\beta-1)}\,dt.
\]
The second integral is finite exactly when $p(\beta-1)>-1$.
In that case,
\[
 u_\beta\in W^{1,p}(I)
 \quad\Longleftrightarrow\quad
 \beta>1-\frac{1}{p},
 \qquad
 \|v_\beta\|_{L^p(I)}^p
 =\frac{2\beta^p}{p(\beta-1)+1}.
\]
At the threshold the exponent is $-1$ and the integral diverges;
below it, the integral also diverges. For $p=\infty$, the derivative
is essentially bounded exactly when $\beta\geq1$.
In particular, $u_1(t)=|t|$ belongs to $W^{1,p}(I)$ for every
$1\leq p\leq\infty$, although it is not $C^1$.
On the other hand, $|t|^{1/4}$ is continuous and does not belong to
$W^{1,2}(I)$. Continuity and membership in a Sobolev space are not
equivalent conditions.
\end{example}

\begin{example}[A jump and the distinction between a class and a representative]
\label{ej:sobolev-salto-y-representantes}
On $I=(-1,1)$, let $H(t):=\mathbf 1_{(0,1)}(t)$.
Although its classical derivative is zero away from $0$, for
$\varphi\in C_c^\infty(I)$ one has
\[
 \int_I H\varphi'=\int_0^1\varphi'=-\varphi(0).
\]
If a weak derivative $v\in L^1_{\mathrm{loc}}(I)$ existed, it would have
to satisfy $\displaystyle\int_Iv\varphi=\varphi(0)$ for every test function.
Choose $\eta\in C_c^\infty((-1,1))$ with $0\leq\eta\leq1$ and $\eta(0)=1$,
and set $\eta_\varepsilon(t):=\eta(t/\varepsilon)$.
Then
\[
 1=\left|\int_Iv\eta_\varepsilon\right|
 \leq\int_{-\varepsilon}^{\varepsilon}|v(t)|\,dt
 \longrightarrow0,
\]
a contradiction. Convergence follows from dominated convergence applied
to $|v|\mathbf 1_{(-\varepsilon,\varepsilon)}$ on an interval compactly
contained in $I$.
Thus $H$ does not belong to $W^{1,p}(I)$ for any $p\in[1,\infty]$.
Changing its value at $0$ affects neither the integrals nor this conclusion.

By contrast, defining $\widetilde u(t)=|t|$ for $t\neq0$ and
$\widetilde u(0)=1$ gives a discontinuous function representing the same
element of $W^{1,p}(I)$ as $|t|$. The spaces $L^p$ and Sobolev spaces
identify functions that agree almost everywhere.
Saying that a Sobolev element is continuous means that it has a continuous
representative, not that every representative is continuous.
The jump of $H$ between two intervals of positive measure cannot be
repaired by modifying it on a null set.
\end{example}

\begin{example}[A Sobolev function with no continuous representative]
\label{ej:sobolev-singularidad-radial}
Let $n\geq2$, $B:=B_{\mathrm{euc}}(0,1)\subseteq\mathbb R^n$, and
$0<\gamma<n-1$. Define $u(x):=|x|^{-\gamma}$ for $x\neq0$, and $u(0):=0$.
Its weak derivative is
\[
 D_i u(x)=-\gamma x_i|x|^{-\gamma-2}
 \quad(x\neq0).
\]
To verify this, integrate by parts on
$B\setminus\overline{B_{\mathrm{euc}}(0,\varepsilon)}$ against
$\varphi\in C_c^\infty(B)$. The outer-sphere term vanishes because
of the support of $\varphi$. The inner term has absolute value at most
\[
 \omega_{n-1}\varepsilon^{n-1-\gamma}\|\varphi\|_\infty,
\]
where $\omega_{n-1}$ is the area of the unit sphere in $\mathbb R^n$.
This bound tends to zero because $\gamma<n-1$. The volume terms converge
because $u$ and $D_i u$ are locally integrable: their radial bounds involve
$\displaystyle\int_0^1 r^{n-1-\gamma}\,dr$ and
$\displaystyle\int_0^1 r^{n-2-\gamma}\,dr$.
Thus the weak-derivative identity persists upon taking the limit.

The radial volume factor can be verified without an angular formula.
If $b_n$ is the volume of the unit ball, dilations give
$\lambda_n(B_{\mathrm{euc}}(0,r))=b_nr^n$.
Thus an annulus with radii $a<b$ has measure
$b_n(b^n-a^n)=nb_n\displaystyle\int_a^b r^{n-1}\,dr$.
Applying this identity first to step functions of $r$ and then,
by monotone convergence, to nonnegative functions gives radial
integration with factor $nb_nr^{n-1}\,dr$.
Moreover, $nb_n=\omega_{n-1}$ by the divergence theorem applied to
the field $x\mapsto x$ on the ball: its divergence is $n$ and its inner
product with the outward unit normal is $1$ on the sphere.
Hence, for $1\leq p<\infty$,
\[
 \int_B|u|^p=\omega_{n-1}\int_0^1r^{n-1-\gamma p}\,dr,
 \qquad
 \int_B|Du|_2^p
 =\omega_{n-1}\gamma^p\int_0^1r^{n-1-p(\gamma+1)}\,dr.
\]
The Euclidean norm $|Du|_2$ and the norm of its components in
$\mathbb R^n$ are equivalent; thus the second integral is finite if
and only if every component $D_i u$ belongs to $L^p(B)$.
Consequently,
\[
 u\in W^{1,p}(B)
 \quad\Longleftrightarrow\quad p(\gamma+1)<n.
\]
This condition also implies $\gamma p<n$, which ensures $u\in L^p$.
When $n=3$, $p=2$, and $\gamma=\displaystyle\frac{1}{4}$,
\[
 \|u\|_{L^2(B)}^2=\frac{8\pi}{5},
 \qquad
 \|Du\|_{L^2(B)}^2=\frac{\pi}{2},
 \qquad
 \|u\|_{W^{1,2}(B)}^2=\frac{21\pi}{10}.
\]
However, no representative is continuous in a neighborhood of $0$.
Every continuous representative would be bounded on a smaller closed
ball, whereas $|x|^{-\gamma}$ exceeds any bound on a sufficiently small
punctured ball of positive measure.
Integrability of the singularity depends both on its strength and on
the dimension.
\end{example}

\begin{figura}[htbp]
 \includegraphics[width=0.98\linewidth]{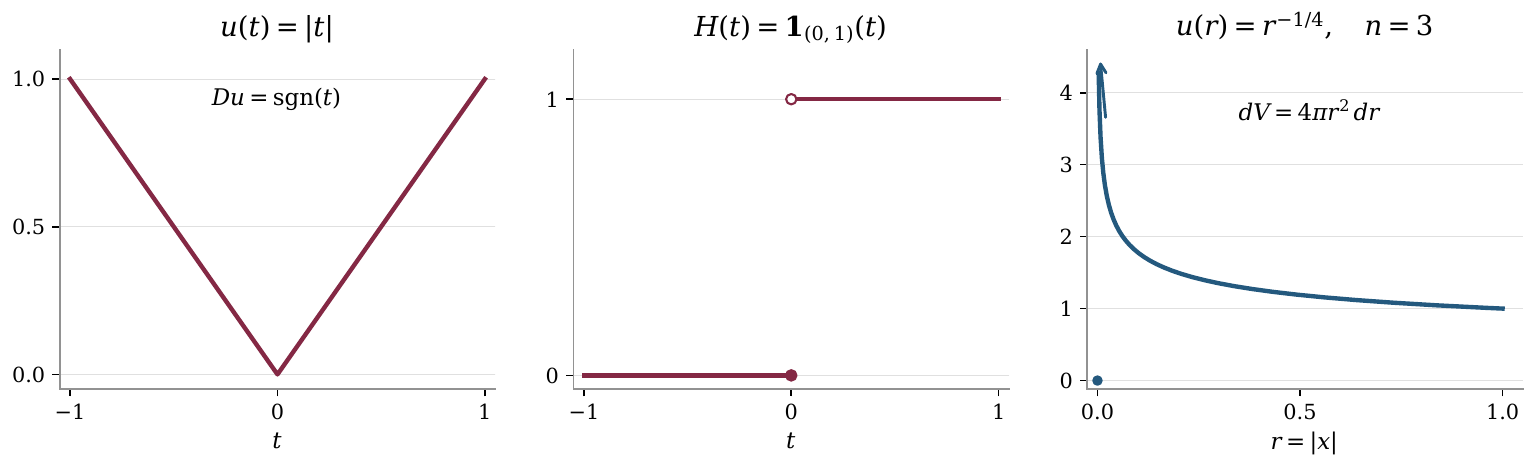}
 \caption{Three different profiles: the corner $|t|$ has a bounded weak
 derivative; the jump $H$ has no locally integrable weak derivative;
 the radial profile $r^{-1/4}$ represents a $W^{1,2}$ function on a ball
 in $\mathbb R^3$. In the third panel, integrability uses the volume
 factor $r^2\,dr$, rather than the measure of an interval.
 The isolated value assigned at the origin does not remove the singularity.}
 \label{fig:perfiles-funciones-sobolev}
\end{figura}

\subsection{Properties of the normed space}

The preceding expression defines a norm on $W^{1,p}(\Omega)$.
\begin{proposition}\label{W1p es normado}
$W^{1,p}(\Omega)$ is a real vector space and $\|\cdot\|_{W^{1,p}(\Omega)}$ is a norm on $W^{1,p}(\Omega)$ for every $p\in [1,\infty]$.
\end{proposition}
\begin{proof}
By Proposition~\ref{linealidad derivada debil}, $W^{1,p}(\Omega)$ is a vector subspace of $L^{p}(\Omega)$, and hence is a vector space. Now, given $u\in W^{1,p}(\Omega)$, if $u=0$, then $D_{i}u=0$ for every $i\in\{1,\dots,n\}$, so that $\|u\|_{W^{1,p}(\Omega)}=0$. Conversely, if $\|u\|_{W^{1,p}(\Omega)}=0$, since $\|u\|_{L^{p}(\Omega)}\leq \|u\|_{W^{1,p}(\Omega)}$, we have $u=0$ in $L^{p}(\Omega)$.

We have $\|\lambda u\|_{W^{1,p}(\Omega)}=|\lambda|\|u\|_{W^{1,p}(\Omega)}$ for any $u\in W^{1,p}(\Omega)$ and $\lambda\in\mathbb{R}$.

Finally, if $p\in [1,\infty)$, applying first Minkowski's inequality in $L^p(\Omega)$, given in Theorem~\ref{teo:b6-espacios-lp-espacios-de-lebesgue}, and then the triangle inequality for the $\|\cdot\|_{p}$ norm on $\mathbb{R}^{n+1}$, we obtain, for any $u,v\in W^{1,p}(\Omega)$,
\[
 \begin{aligned}
 \|u+v\|_{W^{1,p}(\Omega)}
 &\leq\left((\|u\|_{L^p(\Omega)}+\|v\|_{L^p(\Omega)})^p
 +\displaystyle\sum_{i=1}^{n}(\|D_i u\|_{L^p(\Omega)}+\|D_i v\|_{L^p(\Omega)})^p\right)^{\frac{1}{p}}\\
 &\leq\left(\|u\|_{L^p(\Omega)}^p+\displaystyle\sum_{i=1}^{n}\|D_i u\|_{L^p(\Omega)}^p\right)^{\frac{1}{p}}
 +\left(\|v\|_{L^p(\Omega)}^p+\displaystyle\sum_{i=1}^{n}\|D_i v\|_{L^p(\Omega)}^p\right)^{\frac{1}{p}}\\
 &=\|u\|_{W^{1,p}(\Omega)}+\|v\|_{W^{1,p}(\Omega)}.
 \end{aligned}
 \]

If $p=\infty$, the triangle inequality follows in exactly the same way from the triangle inequality for the $\|\cdot\|_{\infty}$ norm on $L^{\infty}(\Omega)$ and that for the $\|\cdot\|_{\infty}$ norm on $\mathbb{R}^{n+1}$.
\end{proof}
The following lemma will be particularly useful for working with Sobolev spaces, both in Euclidean spaces and later on Riemannian manifolds.
\begin{lemma}\label{lema 16.13}
Let $p\in[1,\infty]$ and let $(u_{k})$ be a sequence in $W^{1,p}(\Omega)$ such that $u_{k}\to u$ in $L^{p}(\Omega)$ and $D_{i}u_{k}\to v_{i}$ in $L^{p}(\Omega)$ for each $i\in \{1,\dots,n\}$. Then $u$ is weakly differentiable on $\Omega$, $v_{i}=D_{i}u$ for every $i\in \{1,\dots,n\}$, and $u_{k}\to u$ in $W^{1,p}(\Omega)$.
\end{lemma}
\begin{proof}
Let $\phi\in C_{c}^{\infty}(\Omega)$. Then $\phi,\displaystyle\frac{\partial \phi}{ \partial x_{i}}\in L^{q}(\Omega)$ for every $q\in [1,\infty]$. Let $q> 1$ be such that $\displaystyle\frac{1}{p}+\displaystyle\frac{1}{q}=1$ if $p\in (1,\infty)$, $q=1$ if $p=\infty$, and $q=\infty$ if $p=1$. By Hölder's inequality (Proposition~\ref{desigualdad de holder}), we have \[\left|\int_{\Omega}u\frac{\partial \phi}{\partial x_{i}}-\int_{\Omega}u_{k}\frac{\partial \phi}{\partial x_{i}}\right|\leq \|u-u_{k}\|_{L^p(\Omega)}\left\|\frac{\partial \phi}{\partial x_{i}}\right\|_{L^q(\Omega)}\to 0\text{ as $k\to\infty$},\]

\[\left|\int_{\Omega}v_{i}\phi-\int_{\Omega}(D_{i}u_{k})\phi\right|\leq \left\|v_{i}-D_{i}u_{k}\right\|_{L^p(\Omega)}\|\phi\|_{L^q(\Omega)}\to0\text{ as $k\to\infty$}.\]
Consequently, \[\int_{\Omega}u\frac{\partial \phi}{\partial x_{i}}+\int_{\Omega}v_{i}\phi=\lim_{k\to\infty}\left(\int_{\Omega}u_{k}\frac{\partial\phi}{\partial x_{i}}+\int_{\Omega}(D_{i}u_{k})\phi\right)=0\] for each $i\in\{1,\dots,n\}$. This proves that $u$ is weakly differentiable on $\Omega$ and that $v_{i}=D_{i}u$. Consequently, $u\in W^{1,p}(\Omega)$. Moreover, if $p\in [1,\infty)$, we have \[\lim_{k\to\infty}\|u_{k}-u\|_{W^{1,p}(\Omega)}^{p}=\lim_{k\to\infty}\|u_{k}-u\|_{L^p(\Omega)}^{p}+\displaystyle\sum_{i=1}^{n}\lim_{k\to\infty}\|D_{i}u_{k}-D_{i}u\|_{L^p(\Omega)}^{p}=0.\]

If $p=\infty$, we have
\[
\begin{aligned}
\lim_{k\to\infty}\|u_k-u\|_{W^{1,\infty}(\Omega)}
&=\lim_{k\to\infty}\max\Big\{\|u_k-u\|_{L^\infty(\Omega)},
\|D_1u_k-D_1u\|_{L^\infty(\Omega)},\dots,
\|D_nu_k-D_nu\|_{L^\infty(\Omega)}\Big\}\\
&=0,
\end{aligned}
\]
by continuity of the function $\displaystyle\max$ and the assumed convergences in $L^\infty(\Omega)$.
\end{proof}
The preceding lemma helps us prove the following theorem.
\begin{theorem}\label{teo:sobolev-en-abiertos-euclidianos-espacio-banach}\index{Sobolev space!completeness of \(W^{1,p}(\Omega)\)}
The space $W^{1,p}(\Omega)$ is a Banach space for every $p\in[1,\infty]$.
\end{theorem}
\begin{proof}
Let $(u_{k})$ be a Cauchy sequence in $W^{1,p}(\Omega)$. Then the sequences $(u_{k})$ and $(D_{i}u_{k})$ are Cauchy in $L^{p}(\Omega)$. Since $L^{p}(\Omega)$ is complete, we have $u_{k}\to u$ in $L^{p}(\Omega)$ and $D_{i}u_{k}\to v_{i}$ in $L^{p}(\Omega)$ for each $i\in \{1,\dots,n\}$. Lemma~\ref{lema 16.13} gives $u\in W^{1,p}(\Omega)$ and $u_{k}\to u$ in $W^{1,p}(\Omega)$. This proves that $W^{1,p}(\Omega)$ is complete; since Proposition~\ref{W1p es normado} already shows that it is a normed space, we conclude that it is a Banach space.
\end{proof}

\begin{example}[A smooth Cauchy sequence with no smooth limit]
\label{ej:sucesion-cauchy-suave-no-completa}
\index{Sobolev space!motivation for completion}
Consider the open set $I=(-1,1)\subseteq\mathbb R$ and fix
$1\leq p<\infty$. Write
$\mathcal S_p(I):=C^\infty(I)\cap W^{1,p}(I)$, equipped with
the Sobolev norm. For its elements,
\[
 \|f\|_{W^{1,p}(I)}^p
 =\int_{-1}^1|f(t)|^p\,dt+\int_{-1}^1|f'(t)|^p\,dt.
\]
For $k\in\mathbb N$ define
\[
 u_k(t):=\sqrt{t^2+k^{-2}},\qquad
 u_k'(t)=\frac{t}{\sqrt{t^2+k^{-2}}}.
\]
The argument of the square root is positive throughout $\mathbb R$,
so $u_k$ is smooth. Moreover, $|u_k|\leq\sqrt2$ and $|u_k'|\leq1$
on $I$; thus $u_k\in\mathcal S_p(I)$.

Set $u(t)=|t|$ and $v(t)=\operatorname{sgn}(t)$, with $v(0)=0$.
The identity
\[
 0\leq u_k(t)-u(t)
 =\frac{k^{-2}}{\sqrt{t^2+k^{-2}}+|t|}\leq k^{-1}
\]
gives $\|u_k-u\|_{L^p(I)}^p\leq2k^{-p}$. To estimate the derivatives,
fix $0<a<1$. On $|t|<a$ we have $|u_k'(t)-v(t)|\leq1$.
On $|t|\geq a$,
\[
 |u_k'(t)-v(t)|
 =1-\frac{|t|}{\sqrt{t^2+k^{-2}}}
 =\frac{k^{-2}}
 {\sqrt{t^2+k^{-2}}(\sqrt{t^2+k^{-2}}+|t|)}
 \leq\frac{k^{-2}}{2a^2}.
\]
The two regions have lengths $2a$ and $2-2a$, respectively.
Hence for $k\geq2$ and $a=k^{-1/2}$,
\[
 \|u_k'-v\|_{L^p(I)}^p
 \leq2a+2\left(\frac{k^{-2}}{2a^2}\right)^p
 =2k^{-1/2}+2(2k)^{-p}.
\]
Define
\[
 b_k:=\left(2k^{-p}+2k^{-1/2}+2(2k)^{-p}\right)^{\frac{1}{p}}.
\]
The triangle inequality in $L^p(I)\times L^p(I)$, with the norm
given by the sum of $p$th powers, yields, for $k,\ell\geq2$,
\[
 \begin{aligned}
 \|u_k-u_\ell\|_{W^{1,p}(I)}
 &\leq
 \left(\|u_k-u\|_{L^p(I)}^p+\|u_k'-v\|_{L^p(I)}^p\right)^{\frac{1}{p}}
 \\&\quad+\left(\|u_\ell-u\|_{L^p(I)}^p+
               \|u_\ell'-v\|_{L^p(I)}^p\right)^{\frac{1}{p}}
 \leq b_k+b_\ell.
 \end{aligned}
\]
Since $b_k\to0$, the sequence is Cauchy in the Sobolev norm.
This conclusion only uses the norms of the smooth functions and
the convergence of their pairs $(u_k,u_k')$.

Suppose it converges in that norm to some $w\in\mathcal S_p(I)$.
The Sobolev norm controls the $L^p$ norm, so $w=u$ almost everywhere.
Both functions are continuous. If they differed at a point, their
difference would have positive absolute value on an open interval,
which has positive measure; this contradicts equality almost
everywhere. Therefore $w(t)=|t|$ for all $t\in I$.
However, its difference quotients at zero are $1$ from the right
and $-1$ from the left. No such smooth function $w$ exists.
We have proved that $\mathcal S_p(I)$, equipped with the Sobolev norm,
is an incomplete normed space.

The element added by completion can be identified concretely.
For $\varphi\in C_c^\infty(I)$, integrating over the two intervals,
\[
 \begin{aligned}
 \int_{-1}^1|t|\varphi'(t)\,dt
 &=\int_{-1}^0(-t)\varphi'(t)\,dt
       +\int_0^1t\varphi'(t)\,dt\\
 &=\int_{-1}^0\varphi(t)\,dt-\int_0^1\varphi(t)\,dt
 =-\int_{-1}^1v(t)\varphi(t)\,dt.
 \end{aligned}
\]
The terms at $\pm1$ vanish because of the test function's support,
and those at zero because of the factor $t$. Thus $u$ has weak
derivative $v$, and
\[
 \|u\|_{W^{1,p}(I)}^p
 =2\int_0^1t^p\,dt+2=\frac{2}{p+1}+2.
\]
The preceding estimates prove $u_k\to u$ in $W^{1,p}(I)$.
In particular, $u$ belongs to the closure of $\mathcal S_p(I)$
in that norm and does not belong to $\mathcal S_p(I)$: this
subspace of smooth functions is not closed in $W^{1,p}(I)$.
\end{example}

In addition to being a Banach space, for $p=2$, $W^{1,2}(\Omega)$ is a Hilbert space, for which we use the notation $H^{1}(\Omega):=W^{1,2}(\Omega)$. We record this as a proposition.
\begin{proposition}\label{prop:sobolev-en-abiertos-euclidianos-espacio-hilbert-producto-interior}
If $p=2$, $W^{1,p}(\Omega)$ is a Hilbert space with the inner product \[\langle u,v\rangle_{H^{1}(\Omega)}:=\int_{\Omega}uv+\displaystyle\sum_{i=1}^{n}\int_{\Omega}(D_{i}u)(D_{i}v)=\int_{\Omega}uv+\int_{\Omega}\nabla u\cdot \nabla v.\]
\end{proposition}
\begin{proof}
Hölder's inequality makes all products in the integrand integrable.
The expression is bilinear and symmetric, and
\[
\langle u,u\rangle_{H^1(\Omega)}
=\|u\|_{L^2(\Omega)}^2
+\sum_{i=1}^n\|D_i u\|_{L^2(\Omega)}^2
=\|u\|_{W^{1,2}(\Omega)}^2.
\]
It is therefore nonnegative and vanishes only when $u=0$ as an equivalence
class. It is an inner product whose norm is exactly the Sobolev norm.
Theorem~\ref{teo:sobolev-en-abiertos-euclidianos-espacio-banach}
states that $W^{1,2}(\Omega)$ is complete for this norm, so the space
with the stated inner product is a Hilbert space.
\end{proof}
\begin{remark}\label{obs:sobolev-en-abiertos-euclidianos-funcion-clase-debilmente-diferenciable-sus-derivadas}
$C_{c}^{\infty}(\Omega)\subseteq W^{1,p}(\Omega)$ for every $p\in [1,\infty)$, since every function of class $C_{c}^{\infty}(\Omega)$ is weakly differentiable and its partial derivatives are continuous and compactly supported, and thus belong to $L^{p}(\Omega)$.
\end{remark}
The space $C_{c}^{\infty}(\Omega)$ leads to a basic question in Sobolev analysis: under what conditions can a Sobolev function be approximated by smooth functions with compact support?
We know that, in general, a Sobolev function need not be classically differentiable; nevertheless, under suitable assumptions on the domain, these functions can be approximated by elements of $C_{c}^{\infty}(\Omega)$.
Questions of this kind are called \textit{approximation problems} and provide an essential tool both in the abstract theory and in applications to partial differential equations. We therefore give the following definition.

\begin{definition}\label{def:sobolev-en-abiertos-euclidianos-espacio-de-sobolev-espacio-sobolev-cerradura-denoram}\index{Sobolev space}
Let $p\in[1,\infty)$. We define the Sobolev space $W_{0}^{1,p}(\Omega)$ as the closure of $C_{c}^{\infty}(\Omega)$ in $W^{1,p}(\Omega)$. We write \[H_{0}^{1}(\Omega):=W_{0}^{1,2}(\Omega).\]
\end{definition}
\begin{remark}\label{obs:sobolev-en-abiertos-euclidianos-notacion-anterior-sentido-subespacio-vectorial-cerra}
This notation makes sense because $W_{0}^{1,2}(\Omega)=\overline{C_{c}^{\infty}(\Omega)}^{W^{1,2}(\Omega)}$ is a closed vector subspace of $W^{1,2}(\Omega)$, which is a Hilbert space, and therefore $H_{0}^{1}(\Omega)$ is also a Hilbert space.
\end{remark}
In the case $\Omega=\mathbb{R}^{n}$, there are no boundary obstructions. The argument combines cutoffs moving out toward infinity with regularization by mollifiers, and will serve as the model for the proof on complete Riemannian manifolds.
\begin{proposition}\label{Cc es denso en W1p}
$C_{c}^{\infty}(\mathbb{R}^{n})$ is dense in $W^{1,p}(\mathbb{R}^{n})$ for every $p\in [1,\infty)$.
\end{proposition}

\begin{proof}
Let $u\in W^{1,p}(\mathbb R^n)$. By Proposition~\ref{funciones flan}, there exists $\chi\in C_c^\infty(\mathbb R^n)$ such that $0\leq\chi\leq1$, $\chi=1$ on $B_{\mathrm{euc}}(0,1)$, and $\operatorname{supp}\chi\subseteq B_{\mathrm{euc}}(0,2)$. For $R\geq1$, set
\[
\chi_R(x):=\chi\left(\frac{x}{R}\right),
\qquad
u_R:=\chi_Ru.
\]
The weak Leibniz rule gives
\[
D_i u_R=\chi_RD_i u+uD_i\chi_R.
\]
Indeed, if $\varphi\in C_c^\infty(\mathbb R^n)$, we apply the definition of the weak derivative of $u$ to $\chi_R\varphi$ and expand $D_i(\chi_R\varphi)$. Thus, we do not assume an unproved product rule.

Since $\chi_R(x)\to1$ for every $x\in\mathbb R^n$ and $|1-\chi_R|\leq1$, the dominated convergence theorem \ref{convergencia dominada}, applied to $|u|^p$ and $|D_i u|^p$, shows that
\[
\|(1-\chi_R)u\|_{L^p(\mathbb R^n)}
+
\|(1-\chi_R)D_i u\|_{L^p(\mathbb R^n)}
\longrightarrow0.
\]
Moreover,
\[
\|uD_i\chi_R\|_{L^p(\mathbb R^n)}
\leq
\frac{\|D_i\chi\|_{L^\infty(\mathbb R^n)}}{R}
\|u\|_{L^p(\mathbb R^n)}
\longrightarrow0.
\]
Therefore, $u_R\to u$ in $W^{1,p}(\mathbb R^n)$.

Now fix $R$. Let $\rho$ be the mollifier from Theorem~\ref{teo:molificadores-euclidianos} and define $u_{R,\varepsilon}:=\rho_\varepsilon*u_R$. Then $u_{R,\varepsilon}\in C_c^\infty(\mathbb R^n)$ and, by Fubini's theorem \ref{Fubini} and the definition of weak derivative,
\[
D_i u_{R,\varepsilon}
=
\rho_\varepsilon*D_i u_R.
\]
Continuity of translations from Proposition~\ref{prop:continuidad-traslaciones-Lp}, Minkowski's integral inequality from Theorem~\ref{teo:minkowski-integral}, and the dominated convergence theorem~\ref{convergencia dominada} imply that, for every $v\in L^p(\mathbb R^n)$,
\[
\|\rho_\varepsilon*v-v\|_{L^p(\mathbb R^n)}
\leq
\int_{\mathbb R^n}\rho(z)
\|v(\,\cdot-\varepsilon z)-v\|_{L^p(\mathbb R^n)}\,dz
\longrightarrow0.
\]
We apply this estimate to $v=u_R$ and $v=D_i u_R$ to conclude that $u_{R,\varepsilon}\to u_R$ in $W^{1,p}(\mathbb R^n)$. Choosing first $R\to\infty$ and then $\varepsilon\to0^+$, or taking a diagonal sequence, completes the proof.
\end{proof}

\section{Sobolev embeddings}
By definition, the Sobolev spaces $W_{0}^{1,p}(\Omega)$ are contained in $L^{p}(\Omega)$. An interesting question is: for which other values $q$ is $W_{0}^{1,p}(\Omega)$ continuously embedded in $L^{q}(\Omega)$?

The fundamental tool for answering this question consists of inequalities that estimate the $L^{q}(\mathbb R^n)$ norm of a function $\phi\in C_{c}^{\infty}(\mathbb{R}^{n})$ in terms of the $(L^{p}(\mathbb R^n))^n$ norm of its gradient. Such inequalities are usually called \textit{Sobolev inequalities}.

We begin by giving a norm on the product space
\[(L^{p}(\mathbb{R}^{n}))^{n}:=L^{p}(\mathbb{R}^{n})\times \dots\times L^{p}(\mathbb{R}^{n})\]
with $p\in[1,\infty)$. This norm is given by
\[
\|(f_{1},\dots,f_{n})\|_{(L^p(\mathbb R^n))^n}
:=\left(\displaystyle\sum_{i=1}^{n}\|f_i\|_{L^p(\mathbb R^n)}^p\right)^{\frac{1}{p}}.
\]
Now consider $\varphi\in C_{c}^{\infty}(\mathbb{R}^{n})$. Since $C_{c}^{\infty}(\mathbb{R}^{n})\subseteq L^{p}(\mathbb{R}^{n})$, \[\nabla\varphi=\left(\frac{\partial\varphi}{\partial x_{1}},\dots,\frac{\partial \varphi}{\partial x_{n}}\right)\in (L^{p}(\mathbb{R}^{n}))^{n}\] for every $p\in[1,\infty)$.

We ask for which numbers $p,q\in [1,\infty)$ it is possible to bound $\|\varphi\|_{L^q(\mathbb R^n)}$ uniformly in terms of $\|\nabla\varphi\|_{(L^p(\mathbb R^n))^n}$. We first show that $p$ and $q$ cannot be arbitrary.
\begin{proposition}\label{prop:sobolev-en-abiertos-euclidianos-constante-depende-unicamente-tal}
Let $p,q\in[1,\infty)$. If there is a constant $C>0$ depending only on $n$, $p$, and $q$ such that
\[
 \|\varphi\|_{L^q(\mathbb R^n)}
 \leq C\|\nabla\varphi\|_{(L^p(\mathbb R^n))^n}
 \]
for every $\varphi\in C_{c}^{\infty}(\mathbb{R}^{n})$, then $p<n$ and $q=\displaystyle\frac{np}{n-p}$.
\end{proposition}
\begin{proof}
Fix a nonzero function $\varphi\in C_c^\infty(\mathbb R^n)$ and,
for $\lambda>0$, define $\varphi_\lambda(x):=\varphi(\lambda x)$.
The change of variables $y=\lambda x$ and the chain rule give
\[
\|\varphi_\lambda\|_{L^q(\mathbb R^n)}
=\lambda^{-\frac{n}{q}}\|\varphi\|_{L^q(\mathbb R^n)},\qquad
\|\nabla\varphi_\lambda\|_{(L^p(\mathbb R^n))^n}
=\lambda^{1-\frac{n}{p}}\|\nabla\varphi\|_{(L^p(\mathbb R^n))^n}.
\]
The second norm is positive: if every derivative vanished, $\varphi$
would be constant on $\mathbb R^n$, and compact support would force it
to be zero. Apply the assumed inequality to $\varphi_\lambda$ and
divide by its scaling factor to obtain
\[
0<\|\varphi\|_{L^q(\mathbb R^n)}
\leq C\lambda^{\beta}
\|\nabla\varphi\|_{(L^p(\mathbb R^n))^n},
\qquad \beta:=1-\frac{n}{p}+\frac{n}{q}.
\]
If $\beta>0$, letting $\lambda\to0^+$ contradicts this inequality;
if $\beta<0$, letting $\lambda\to\infty$ gives the same contradiction.
Thus $\beta=0$ and
\[
\frac{n}{p}=1+\frac{n}{q},\qquad q=\frac{np}{n-p}.
\]
Since $q<\infty$, the first equality implies $\displaystyle\frac{n}{p}>1$, or $p<n$.
\end{proof}
This value of $q$ plays an important role in Sobolev space theory and thus receives a special name.
\begin{definition}\label{def:sobolev-en-abiertos-euclidianos-define-exponente-critico-sobolev-notemos}\index{critical Sobolev exponent}
If $p\in [1,n)$, the critical Sobolev exponent is defined by \[p^{*}:=\frac{np}{n-p}.\]
Moreover, $p^{*}>p$.
\end{definition}
In fact, the inequality involving the norm of the gradient holds for $p\in [1,n)$ and $q=p^{*}$, as we shall now prove. We first need a lemma.
\begin{lemma}\label{lema 5 para gagliardo}
Let $n\geq 2$ and $f_{1},\dots, f_{n}\in L^{n-1}(\mathbb{R}^{n-1})$. If $x=(x_{1},\dots,x_{n})\in\mathbb{R}^{n}$, write $\hat{x}_{i}:=(x_{1},\dots,x_{i-1},x_{i+1},\dots,x_{n})\in \mathbb{R}^{n-1}$ and define \[f(x):=\prod_{i=1}^{n}f_{i}(\hat{x}_{i})\] Then $f\in L^{1}(\mathbb{R}^{n})$ and \[\|f\|_{L^{1}(\mathbb{R}^{n})}\leq \prod_{i=1}^{n}\|f_{i}\|_{L^{n-1}(\mathbb{R}^{n-1})}.\]
\end{lemma}
\begin{proof}
We argue by induction on $n$. If $n=2$, we have \[\int_{\mathbb{R}^{2}}|f(x_{1},x_{2})|dx_{1}dx_{2}=\int_{\mathbb{R}^{2}}|f_{1}(x_{2})f_{2}(x_{1})|dx_{1}dx_{2}=\int_{\mathbb{R}^{2}}|f_{1}(x_{2})||f_{2}(x_{1})|dx_{1}dx_{2}\] \[=\int_{\mathbb{R}}|f_{1}(x_{2})|dx_{2}\int_{\mathbb{R}}|f_{2}(x_{1})|dx_{1}.\] Suppose that the assertion holds for some $n\geq2$; we show that it holds for $n+1$.

Let $f_{1},\dots,f_{n+1}\in L^{n}(\mathbb{R}^{n})$. Fix the value of $x_{n+1}$ and define \[g_{i}(z_{1},\dots,z_{n-1}):=|f_{i}(z_{1},\dots,z_{n-1},x_{n+1})|^{\frac{n}{n-1}},i\in\{1,\dots,n\}.\]

By Theorem~\ref{Fubini}, we have $g_{i}\in L^{n-1}(\mathbb{R}^{n-1})$ for almost every $x_{n+1}$ (by Fubini's theorem with $m=n-1$). The induction hypothesis therefore shows that the function \[g(y):=\prod_{i=1}^{n}g_{i}(\hat{y}_{i})=\prod_{i=1}^{n}|f_{i}(\hat{y}_{i},x_{n+1})|^{\frac{n}{n-1}}=\left(\prod_{i=1}^{n}|f_{i}(\hat{y}_{i},x_{n+1})|\right)^{\frac{n}{n-1}}\] $y\in \mathbb{R}^{n}$ belongs to $L^{1}(\mathbb{R}^{n})$ and that \[\|g\|_{L^{1}(\mathbb{R}^{n})}\leq \prod_{i=1}^{n}\|g_{i}\|_{L^{n-1}(\mathbb{R}^{n-1})}=\prod_{i=1}^{n}\left(\int_{\mathbb{R}^{n-1}}|f_{i}(z,x_{n+1})|^{n}dz\right)^{\frac{1}{n-1}}.\]
We have \[|f(y,x_{n+1})|=\left(\prod_{i=1}^{n}|f_{i}(\hat{y}_{i},x_{n+1})|\right)|f_{n+1}(y)|=|g(y)|^{\frac{n-1}{n}}|f_{n+1}(y)|.\]

Using Proposition~\ref{desigualdad de holder} (Hölder's inequality) and the preceding inequality for $\|g\|_{L^{1}(\mathbb{R}^{n})}$, we obtain the following:

\[\int_{\mathbb{R}^{n}}|f(y,x_{n+1})|dy\leq \|g\|_{L^{1}(\mathbb{R}^{n})}^{\frac{n-1}{n}}\|f_{n+1}\|_{L^{n}(\mathbb{R}^{n})}\leq \prod_{i=1}^{n}\left(\int_{\mathbb{R}^{n-1}}|f_{i}(z,x_{n+1})|^{n}dz\right)^{\frac{1}{n}}\|f_{n+1}\|_{L^{n}(\mathbb{R}^{n})}.\]

The preceding argument holds for almost every $x_{n+1}$, so the functions
\[
h_{i}(x_{n+1})
:=
\left(
\int_{\mathbb{R}^{n-1}}
|f_{i}(z,x_{n+1})|^{n}\,dz
\right)^{\frac{1}{n}},
\qquad i\in\{1,\dots,n\},
\]
belong to $L^{n}(\mathbb{R})$. Indeed, by Tonelli's theorem,
\[
\begin{aligned}
\|h_{i}\|_{L^{n}(\mathbb{R})}^{n}
&=
\int_{\mathbb{R}}
\int_{\mathbb{R}^{n-1}}
|f_{i}(z,x_{n+1})|^{n}\,dz\,dx_{n+1}\\
&=\|f_{i}\|_{L^{n}(\mathbb{R}^{n})}^{n}.
\end{aligned}
\]
Since each $h_{i}\in L^{n}(\mathbb{R})$, we may apply Proposition~\ref{desigualdad de holder} with $n$ factors and all exponents equal to $n$. Thus,
\[
\int_{\mathbb{R}}\prod_{i=1}^{n}|h_{i}(x_{n+1})|\,dx_{n+1}
\leq \prod_{i=1}^{n}\|h_{i}\|_{L^{n}(\mathbb{R})}.
\]

We previously showed that \[\int_{\mathbb{R}^{n}}|f(y,x_{n+1})|dy\leq\prod_{i=1}^{n}\left(\int_{\mathbb{R}^{n-1}}|f_{i}(z,x_{n+1})|^{n}dz\right)^{\frac{1}{n}}\|f_{n+1}\|_{L^{n}(\mathbb{R}^{n})}.\] Integrating this inequality with respect to $x_{n+1}$ and using the preceding inequality gives
\[
\begin{aligned}
\int_{\mathbb{R}^{n+1}}|f(x)|\,dx
&\leq
\left(
\int_{\mathbb{R}}
\prod_{i=1}^{n}|h_{i}(x_{n+1})|\,dx_{n+1}
\right)
\|f_{n+1}\|_{L^{n}(\mathbb{R}^{n})}\\
&\leq
\left(\prod_{i=1}^{n}\|h_{i}\|_{L^{n}(\mathbb{R})}\right)
\|f_{n+1}\|_{L^{n}(\mathbb{R}^{n})}\\
&=
\prod_{i=1}^{n+1}\|f_{i}\|_{L^{n}(\mathbb{R}^{n})},
\end{aligned}
\]
which is the desired conclusion.
\end{proof}
\begin{theorem}\label{gagliardo nirenberg sobolev}\index{Gagliardo Nirenberg Sobolev inequality@Gagliardo--Nirenberg--Sobolev inequality}
There are constants $C>0$, depending only on $n$ and $p$, with the following properties: \begin{enumerate}[label=(\alph*)]
 \item Gagliardo--Nirenberg--Sobolev inequality: If $p\in[1,n)$, then \[\|\varphi\|_{L^{p^{*}}(\mathbb R^n)}\leq C\|\nabla\varphi\|_{(L^p(\mathbb R^n))^n}\] for every $\varphi\in C_{c}^{\infty}(\mathbb{R}^{n})$.
 \item If $p\in(n,\infty)$, then \[\|\varphi\|_{L^{\infty}(\mathbb R^n)}\leq C(\lambda_{n}(\supp (\varphi)))^{\frac{1}{n}-\frac{1}{p}}\|\nabla\varphi\|_{(L^p(\mathbb R^n))^n}\] for every $\varphi\in C_{c}^{\infty}(\mathbb{R}^{n})$, where $\lambda_{n}(\supp (\varphi))$ denotes the Lebesgue measure of the support of $\varphi$.
 \end{enumerate}
\end{theorem}
\begin{proof}
First suppose $n\geq2$ and set $r=\displaystyle\frac{n}{n-1}$.
For $\varphi\in C_c^1(\mathbb R^n)$, the fundamental theorem of
calculus in each direction gives
\[
|\varphi(x)|\leq
\int_{\mathbb R}|D_i\varphi(x_1,\ldots,x_{i-1},t,x_{i+1},\ldots,x_n)|\,dt.
\]
Let $\widehat x_i$ omit coordinate $i$ and define
\[
f_i(\widehat x_i):=
\left(\int_{\mathbb R}|D_i\varphi(x_1,\ldots,t,\ldots,x_n)|\,dt\right)^{\frac{1}{n-1}}.
\]
Tonelli's theorem gives $\|f_i\|_{n-1}=\|D_i\varphi\|_1^{\displaystyle\frac{1}{n-1}}$.
Raise each bound to $\displaystyle\frac{1}{n-1}$ and multiply to obtain
$|\varphi|^r\leq\prod_i f_i(\widehat x_i)$.
Integrate and apply Lemma~\ref{lema 5 para gagliardo}.
The geometric--arithmetic mean inequality yields
\[
\|\varphi\|_{L^r(\mathbb R^n)}
\leq\prod_{i=1}^n\|D_i\varphi\|_{L^1(\mathbb R^n)}^{\frac{1}{n}}
\leq\frac{1}{n}\sum_{i=1}^n\|D_i\varphi\|_{L^1(\mathbb R^n)}.
\]
This proves (a) for $p=1$.

For $\delta>1$, the function $h_\delta(t):=t|t|^{\delta-1}$
is $C^1$ on the whole line, with
$h_\delta'(t)=\delta|t|^{\delta-1}$, including at $t=0$.
Hence $h_\delta\circ\varphi\in C_c^1(\mathbb R^n)$.
Apply the preceding bound to this composition, followed by Hölder's
inequality in space and then in the finite sum. For $1<p<\infty$ this gives
\[
\|\varphi\|_{L^{\delta r}(\mathbb R^n)}^\delta
\leq\delta
\left(\int_{\mathbb R^n}|\varphi|^{(\delta-1)p'}\right)^{\frac{1}{p'}}
\|\nabla\varphi\|_{(L^p(\mathbb R^n))^n}.
\]
If $1<p<n$, choose
\[
\delta:=\frac{p(n-1)}{n-p}>1,
\qquad \delta r=(\delta-1)p'=p^*.
\]
For $\varphi\neq0$, divide by
$\|\varphi\|_{L^{p^*}(\mathbb R^n)}^{\delta-1}$ to obtain (a).
For $\varphi=0$, the conclusion is immediate.

Now let $p>n$. For a nonzero function write
\[
S:=|\operatorname{supp}\varphi|>0,
\qquad G:=\|\nabla\varphi\|_{(L^p(\mathbb R^n))^n}>0.
\]
The positivity of $G$ follows because a compactly supported function
with zero gradient is zero. Define
\[
\lambda:=S^{\frac{1}{n}},\qquad
w(x):=\frac{\varphi(\lambda x)}{\lambda^{1-\frac{n}{p}}G}.
\]
Then $|\operatorname{supp}w|=1$ and
$\|\nabla w\|_{(L^p(\mathbb R^n))^n}=1$. The first bound and
Hölder's inequality on this support give $\|w\|_{L^r(\mathbb R^n)}\leq1$.
Set
\[
\gamma:=\frac{r}{p'}>1,\qquad q_k:=r\gamma^k,
\qquad a_k:=\|w\|_{L^{q_k}(\mathbb R^n)},\quad k\in\mathbb N_0.
\]
For $k\geq1$, use the composition with $\delta=\gamma^k$.
The number $t_k:=(\gamma^k-1)p'$ satisfies
$0<t_k=q_{k-1}-p'<q_{k-1}$. Hölder's inequality on a set of measure one
gives directly
\[
\int_{\mathbb R^n}|w|^{t_k}
\leq\left(\int_{\mathbb R^n}|w|^{q_{k-1}}\right)^{\frac{t_k}{q_{k-1}}}.
\]
This inequality also holds when $t_k<1$, since the Hölder exponent
used is $\displaystyle\frac{q_{k-1}}{t_k}>1$. Consequently,
\[
a_k\leq\gamma^{\frac{k}{\gamma^k}}a_{k-1}^{1-\frac{1}{\gamma^k}}.
\]
With $b_k:=\max\{1,a_k\}$, we have $b_0=1$ and
\[
b_k\leq\gamma^{\frac{k}{\gamma^k}}b_{k-1}
\leq\gamma^{\sum_{i=1}^k \frac{i}{\gamma^i}}
\leq\gamma^{\frac{\gamma}{(\gamma-1)^2}}.
\]
Since $q_k\to\infty$, these bounds give the same bound for
$\|w\|_\infty$. Indeed, for any $a<\|w\|_\infty$, the set
$\{|w|>a\}$ has positive measure and
$\|w\|_{q_k}\geq a|\{|w|>a\}|^{\displaystyle\frac{1}{q_k}}$; the upper bound
$\|w\|_{q_k}\leq\|w\|_\infty$ uses the support of measure one.
Undoing the normalization gives
\[
\|\varphi\|_{L^\infty(\mathbb R^n)}
\leq C_{n,p}S^{\frac{1}{n}-\frac{1}{p}}
\|\nabla\varphi\|_{(L^p(\mathbb R^n))^n}.
\]
The zero function also satisfies this bound.

Finally, if $n=1$ and $p>1$, the fundamental theorem of calculus
and Hölder's inequality on the support give
\[
\|\varphi\|_\infty\leq\|\varphi'\|_1
\leq |\operatorname{supp}\varphi|^{1-\frac{1}{p}}\|\varphi'\|_p.
\]
In dimension one this proves (b), and (a) has no admissible indices.
\end{proof}
For $u\in W_{0}^{1,p}(\Omega)$, we shall write $\nabla u=(D_1u,\dots,D_nu)$ and
\[
\|\nabla u\|_{(L^p(\Omega))^n}
:=
\left(\displaystyle\sum_{i=1}^{n}\|D_i u\|_{L^p(\Omega)}^p\right)^{\frac{1}{p}}.
\]
The extension of the preceding inequality to the completed space is recorded in the following corollary.
\begin{corollary}\label{corolario 9}
If $\Omega\subseteq \mathbb{R}^{n}$ is open and $p\in[1,n)$, then $W_{0}^{1,p}(\Omega)\subseteq L^{p^{*}}(\Omega)$, and there is a constant $C>0$, depending only on $n$ and $p$, such that \[\|u\|_{L^{p^{*}}(\Omega)}\leq C\|\nabla u\|_{(L^p(\Omega))^n}\] for every $u\in W_{0}^{1,p}(\Omega)$.
\end{corollary}
\begin{proof}
The Gagliardo--Nirenberg--Sobolev inequality shows that the inclusion operator
\[
I_0\colon C_c^{\infty}(\Omega)\longrightarrow L^{p^*}(\Omega)
\]
is linear and satisfies $\|I_0\varphi\|_{L^{p^*}(\Omega)}\leq C\|\varphi\|_{W^{1,p}(\Omega)}$. Since $W_0^{1,p}(\Omega)$ is the closure of $C_c^{\infty}(\Omega)$ in $W^{1,p}(\Omega)$ and $L^{p^*}(\Omega)$ is a Banach space, Theorem~\ref{teo:extension-operadores-subespacio-denso} provides a unique continuous extension $I\colon W_0^{1,p}(\Omega)\longrightarrow L^{p^*}(\Omega)$.

It remains to identify $Iu$ with the usual representative of $u$ and retain the estimate in terms of the gradient. Let $(\varphi_k)\subseteq C_c^{\infty}(\Omega)$ be a sequence converging to $u$ in $W^{1,p}(\Omega)$. Then $\varphi_k\to u$ in $L^p(\Omega)$ and $\varphi_k\to Iu$ in $L^{p^*}(\Omega)$. Proposition~\ref{Teorema 14.28}, applied to subsequences, shows that both limits agree almost everywhere. Moreover,
\[
\|Iu\|_{L^{p^*}(\Omega)}
=
\lim_{k\to\infty}\|\varphi_k\|_{L^{p^*}(\Omega)}
\leq
C\lim_{k\to\infty}\|\nabla\varphi_k\|_{(L^p(\Omega))^n}
=
C\|\nabla u\|_{(L^p(\Omega))^n}.
\]
Consequently, $u\in L^{p^*}(\Omega)$, the extension agrees with the canonical inclusion, and we obtain the stated inequality.
\end{proof}
The preceding corollary, combined with interpolation between $L^p(\mathbb R^n)$ and $L^{p^*}(\mathbb R^n)$, yields the first continuous Sobolev embedding of this chapter.
\begin{theorem}[Sobolev embedding]\label{encaje de sobolev}\index{Sobolev embedding theorem}
If $\Omega$ is an open subset of $\mathbb{R}^{n}$ and $p\in[1,n)$, then \[W_{0}^{1,p}(\Omega)\subseteq L^{q}(\Omega)\] for every $q\in[p,p^{*}]$, and this inclusion is continuous.
\end{theorem}
\begin{proof}
Corollary~\ref{corolario 9} states that $W_{0}^{1,p}(\Omega)\subseteq L^{p^{*}}(\Omega)$ and that there is a constant $C>0$ such that \[\|u\|_{L^{p^{*}}(\Omega)}\leq C\|\nabla u\|_{(L^p(\Omega))^n}\leq C\|u\|_{W^{1,p}(\Omega)}\] for every $u\in W_{0}^{1,p}(\Omega)$.

Since we also have $W_{0}^{1,p}(\Omega)\subseteq L^{p}(\Omega)$ and $\|u\|_{L^p(\Omega)}\leq \|u\|_{W^{1,p}(\Omega)}$, Theorem~\ref{Ejercicio 14.67} ensures that $u\in L^{q}(\Omega)$ for every $u\in W_{0}^{1,p}(\Omega)$ and $q\in [p,p^{*}]$, and that \[\|u\|_{L^q(\Omega)}\leq\|u\|_{L^p(\Omega)}^{1-\alpha}\|u\|_{L^{p^{*}}(\Omega)}^{\alpha}\leq C^{\alpha}\|u\|_{W^{1,p}(\Omega)}^{1-\alpha}\|u\|_{W^{1,p}(\Omega)}^{\alpha}=C^{\alpha}\|u\|_{W^{1,p}(\Omega)}\], where $\alpha\in[0,1]$ satisfies $\displaystyle\frac{1}{q}=\displaystyle\frac{1-\alpha}{p}+\displaystyle\frac{\alpha}{p^{*}}$. This proves that $W_{0}^{1,p}(\Omega)\subseteq L^{q}(\Omega)$ and that this inclusion is continuous.
\end{proof}
When $\Omega$ is bounded, $\|u\|_{L^q(\Omega)}$ can be bounded in terms of $\|\nabla u\|_{(L^p(\Omega))^n}$ and the measure of $\Omega$ for a wider range of values of $p$ and $q$. The following inequality is known as the Poincaré inequality.
\begin{theorem}[Poincaré inequality]\label{desigualdad de poincaré}\index{Poincare inequality@Poincaré inequality}\glsadd{desigualdad-poincare}
Let $\Omega$ be a bounded open subset of $\mathbb{R}^{n}$. There is a constant $C>0$ depending only on $n$, $p$, and $q$ such that \[\|u\|_{L^q(\Omega)}\leq C(\lambda_{n}(\Omega))^{\frac{1}{q}+\frac{1}{n}-\frac{1}{p}}\|\nabla u\|_{(L^p(\Omega))^n}\] for every $u\in W_{0}^{1,p}(\Omega)$ if any of the following three conditions holds:
 \begin{enumerate}[label=(\alph*)]
 \item $p\in [1,n)$ and $q\in [1,p^{*}]$.
 \item $p=n$ and $q\in [1,\infty)$.
 \item $p\in (n,\infty)$ and $q\in [1,\infty]$.
 \end{enumerate}
Consequently:
 \begin{itemize}
 \item If $p\in [1,n)$, then $W_{0}^{1,p}(\Omega)\subseteq L^{q}(\Omega)$ for each $q\in [1,p^{*}]$, and the inclusion is continuous.
 \item $W_{0}^{1,n}(\Omega)\subseteq L^{q}(\Omega)$ for each $q\in[1,\infty)$, and the inclusion is continuous.
 \item Furthermore, if $p\in(n,\infty)$, then, after choosing a representative, $W_{0}^{1,p}(\Omega)\subseteq C^{0}(\overline{\Omega})$, and the inclusion is continuous.
 \end{itemize}

\end{theorem}
\begin{proof}
\begin{enumerate}[label=(\alph*)]
 \item If $p\in[1,n)$ and $q\in[1,p^{*}]$, applying the inequality from Corollary~\ref{corolario 9} and Proposition~\ref{Proposición 14.31} gives \[\|u\|_{L^q(\Omega)}\leq (\lambda_{n}(\Omega))^{\frac{1}{q}-\frac{1}{p^{*}}}\|u\|_{L^{p^{*}}(\Omega)}\leq C(\lambda_{n}(\Omega))^{\frac{1}{q}+\frac{1}{n}-\frac{1}{p}}\|\nabla u\|_{(L^p(\Omega))^n}\] for every $u\in W_{0}^{1,p}(\Omega)$.
 \item If $n=1$, then for $\varphi\in C_c^\infty(\Omega)$, extension by zero and the fundamental theorem of calculus give $\|\varphi\|_{L^\infty(\Omega)}
\leq\|\varphi'\|_{L^1(\Omega)}$. By Hölder's inequality,
\[
\|\varphi\|_{L^q(\Omega)}
\leq\lambda_1(\Omega)^{\frac{1}{q}}
\|\varphi'\|_{L^1(\Omega)}.
\]
Approximation by test functions extends this estimate to $W_0^{1,1}(\Omega)$ and proves this part in dimension one. Now assume $n\geq2$.
Let $p=n$ and $q\in[1,\infty)$. Define $r:=\displaystyle\max\{\displaystyle\frac{nq}{n+q},1\}$. The number $r$ belongs to $[1,n)$ because $n^{2}>0$ implies $n^{2}+nq>nq$, and hence $\displaystyle\frac{nq}{n+q}<n$. Moreover, $r^{*}=q$ if $r=\displaystyle\frac{nq}{n+q}$, and $q\leq \displaystyle\frac{n}{n-1}$ if $r=1$. Thus, by $(a)$, Proposition~\ref{Proposición 14.31}, and Proposition~\ref{ejercicio 14.70}, for every $u\in W_{0}^{1,p}(\Omega)$ we have \[\|u\|_{L^q(\Omega)}\leq C\|\nabla u\|_{(L^r(\Omega))^n}\leq C n^{\frac{1}{q}}(\lambda_{n}(\Omega))^{\frac{1}{q}}\|\nabla u\|_{(L^n(\Omega))^n}\] if $r=\displaystyle\frac{nq}{n+q}$ and \[\|u\|_{L^q(\Omega)}\leq (\lambda_n(\Omega))^{\frac{1}{q}-\frac{n-1}{n}}\|u\|_{\frac{n}{n-1}}\leq C(\lambda_{n}(\Omega))^{\frac{1}{q}-\frac{n-1}{n}}\|\nabla u\|_{(L^1(\Omega))^n}\leq Cn^{\frac{n-1}{n}}(\lambda_{n}(\Omega))^{\frac{1}{q}}\|\nabla u\|_{(L^n(\Omega))^n}\] if $r=1$.
 \item Let $p\in (n,\infty)$ and $q\in [1,\infty]$. The second inequality in Theorem~\ref{gagliardo nirenberg sobolev}, applied to the extension by zero, gives \[\|\varphi\|_{L^{\infty}(\Omega)}\leq C(\lambda_{n}(\Omega))^{\frac{1}{n}-\frac{1}{p}}\|\nabla\varphi\|_{(L^p(\Omega))^n}\] for every $\varphi\in C_{c}^{\infty}(\Omega)$.

Now recall that $(C^{0}(\overline{\Omega}),\|\cdot\|_{C^0(\overline\Omega)})$ is a Banach space, since $\Omega$ is bounded and therefore $\overline{\Omega}$ is compact. Arguing as in Corollary~\ref{corolario 9} shows that each $u\in W_{0}^{1,p}(\Omega)$ agrees almost everywhere with a function belonging to $C^{0}(\overline{\Omega})$, and that the preceding inequality holds for every $u$. Combining that inequality with Proposition~\ref{Proposición 14.31} gives \[\|u\|_{L^q(\Omega)}\leq (\lambda_{n}(\Omega))^{\frac{1}{q}}\|u\|_{C^0(\overline\Omega)}\leq C(\lambda_{n}(\Omega))^{\frac{1}{q}+\frac{1}{n}-\frac{1}{p}}\|\nabla u\|_{(L^p(\Omega))^n}\], which is the desired inequality.
\end{enumerate}

\end{proof}
Just as we defined the Sobolev spaces $W^{1,p}(\Omega)$, we can define the ``higher-order Sobolev spaces'' $W^{m,p}(\Omega)$:
 \begin{definition}\label{def sobolev generalizado}\index{Sobolev space!higher-order}
Let $\Omega\subseteq \mathbb{R}^{n}$ be open, and let $m\in\mathbb{N}$, $m\geq 2$, and $p\in [1,\infty]$. Define recursively \[W^{m,p}(\Omega):=\{u\in L^{p}(\Omega)\mid u\in W^{1,p}(\Omega)\text{ and }D_{i}u\in W^{m-1,p}(\Omega),\forall i\in\{1,\dots,n\}\}.\]

That is, $u\in W^{m,p}(\Omega)$ if its weak derivatives of order at most $m$, namely \[D^{\alpha}u:=D_{1}^{\alpha_{1}}\cdots D_{n}^{\alpha_{n}}u\] with $\alpha=(\alpha_{1},\dots,\alpha_{n})\in (\mathbb{N}\cup \{0\})^{n}$, $|\alpha|:=\alpha_{1}+\cdots+\alpha_{n}\leq m$, exist and belong to $L^{p}(\Omega)$, where $D_{i}^{0}u:=u$ and $D_{i}^{\alpha_{i}}u:=\underbrace{D_{i}\cdots D_{i}}_{\alpha_{i}-\text{times}}u$ if $\alpha_{i}\geq 1$. Note that all weak derivatives of order at most $m$ have this form.

Using this notation, define
\[
 \|u\|_{W^{m,p}(\Omega)}
 =\left(\displaystyle\sum_{|\alpha|\leq m}\|D^{\alpha}u\|_{L^p(\Omega)}^{p}\right)^{\frac{1}{p}},
 \quad \text{if $p<\infty$},
 \]
and
\[
 \|u\|_{W^{m,\infty}(\Omega)}
 =\max_{|\alpha|\leq m}\|D^{\alpha}u\|_{L^\infty(\Omega)}.
 \]
Finally, define $W_{0}^{m,p}(\Omega)=\overline{C_{c}^{\infty}(\Omega)}^{W^{m,p}(\Omega)}$.
 \end{definition}

The multi-index description agrees with the recursive definition. Indeed, if $D_{i_1}\cdots D_{i_r}u$ exists as an iterated weak derivative, successive applications of the definition give, for $\varphi\in C_c^\infty(\Omega)$,
\[
\int_\Omega D_{i_1}\cdots D_{i_r}u\,\varphi
=(-1)^r\int_\Omega u\,\partial_{i_r}\cdots\partial_{i_1}\varphi.
\]
The classical derivatives of $\varphi$ commute. Uniqueness of weak derivatives therefore implies that two existing iterations with the same indices in different orders agree almost everywhere. Conversely, if all derivatives $D^\alpha u$, with $|\alpha|\leq m$, exist and belong to $L^p(\Omega)$, then for $|\alpha|<m$,
\[
\int_\Omega D^\alpha u\,\partial_i\varphi
=(-1)^{|\alpha|}\int_\Omega u\,\partial^\alpha\partial_i\varphi
=-\int_\Omega D^{\alpha+e_i}u\,\varphi.
\]
Thus $D_i(D^\alpha u)=D^{\alpha+e_i}u$ weakly. Induction on the order supplies every iteration required by the recursive definition.
We have that ($W^{m,p}(\Omega),\|\cdot\|_{W^{m,p
}(\Omega)})$ is a real normed vector space; the proof is analogous to that of Proposition~\ref{W1p es normado}. We shall give results analogous to those for the spaces $W^{1,p}(\Omega)$ and use the latter to establish the corresponding proofs.
 \begin{lemma}\label{lema 16.13 generalizado}
Let $p\in[1,\infty]$ and $m\geq 1$, and let $(u_{k})$ be a sequence in $W^{m,p}(\Omega)$ such that $u_{k}\to u$ in $L^{p}(\Omega)$ and $D^{\alpha}u_{k}\to v_{\alpha}$ in $L^{p}(\Omega)$ for each multi-index $\alpha$ with $|\alpha|\leq m$. Then $u\in W^{m,p}(\Omega)$, $v_{\alpha}=D^{\alpha}u$ for every multi-index $\alpha$ with $|\alpha|\leq m$, and $u_{k}\to u$ in $W^{m,p}(\Omega)$.
\end{lemma}
 \begin{proof}
We prove the assertion by induction on $m$. The case $m=1$ is Lemma~\ref{lema 16.13}. Assume $m\geq2$ and fix a multi-index $\beta$ with $|\beta|\leq m-1$. For each $i\in\{1,\dots,n\}$,
\[
D^\beta u_k\longrightarrow v_\beta,
\qquad
D_i(D^\beta u_k)=D^{\beta+e_i}u_k\longrightarrow v_{\beta+e_i}
\quad\text{on }L^p(\Omega).
\]
Lemma~\ref{lema 16.13}, applied to the sequence $(D^\beta u_k)$, shows that $v_\beta\in W^{1,p}(\Omega)$ and $D_i v_\beta=v_{\beta+e_i}$. For $\beta=0$, uniqueness of the limit in $L^p(\Omega)$ gives $v_0=u$. Iterating these identities according to the length of $\beta$ yields $v_\alpha=D^\alpha u$ for every $|\alpha|\leq m$; in particular, $u\in W^{m,p}(\Omega)$.

If $p<\infty$, then
\[
\|u_k-u\|_{W^{m,p}(\Omega)}^p
=\displaystyle\sum_{|\alpha|\leq m}
\|D^\alpha u_k-D^\alpha u\|_{L^p(\Omega)}^p
\longrightarrow0.
\]
If $p=\infty$, the family of multi-indices with $|\alpha|\leq m$ is finite, and therefore
\[
\|u_k-u\|_{W^{m,\infty}(\Omega)}
=\max_{|\alpha|\leq m}
\|D^\alpha u_k-D^\alpha u\|_{L^\infty(\Omega)}
\longrightarrow0.
\]
This completes the proof.

 \end{proof}

\begin{theorem}\label{teo:sobolev-en-abiertos-euclidianos-espacio-banach-2}\index{Sobolev space!completeness of \(W^{m,p}(\Omega)\)}
$W^{m,p}(\Omega)$ is a Banach space for every $p\in [1,\infty]$.
\end{theorem}
\begin{proof}
Let $(u_{k})$ be a Cauchy sequence in $W^{m,p}(\Omega)$. Then the sequences $(u_{k})$ and $(D^{\alpha}u_{k})$ are Cauchy in $L^{p}(\Omega)$ for every multi-index $\alpha$ with $|\alpha|\leq m$. Since $L^{p}(\Omega)$ is a Banach space by Theorem~\ref{teo:b6-espacios-lp-espacios-de-lebesgue}, we have $u_{k}\to u$ in $L^{p}(\Omega)$ and $D^{\alpha}u_{k}\to v_{\alpha}$ in $L^{p}(\Omega)$ for each multi-index $\alpha$ with $|\alpha|\leq m$. By Lemma~\ref{lema 16.13 generalizado}, $u\in W^{m,p}(\Omega)$ and $u_{k}\to u$ in $W^{m,p}(\Omega)$. Thus $W^{m,p}(\Omega)$ is complete.
\end{proof}

\subsection{Density and extension by zero}

Approximation by cutoffs and mollifiers also controls all derivatives of integer order. The following local construction is the Meyers--Serrin theorem; see \cite[Section~11.2]{Leoni2017}.

\begin{lemma}[Restriction and extension by zero with interior support]
\label{lem:restriccion-extension-local-euclidiana}
\index{extension by zero!with interior support}
Let $U\subseteq V\subseteq\mathbb R^n$ be open, and let $m\in\mathbb N_0$ and $1\leq p<\infty$. The restriction of $w\in W^{m,p}(V)$ belongs to $W^{m,p}(U)$ and satisfies
\[
D^\alpha(w\restriction_U)=(D^\alpha w)\restriction_U,
\quad |\alpha|\leq m,
\qquad
\|w\restriction_U\|_{W^{m,p}(U)}\leq\|w\|_{W^{m,p}(V)}.
\]
If $u\in W^{m,p}(U)$ vanishes almost everywhere outside a compact set $K\subset U$, its extension by zero $\widetilde u$ to $V$ belongs to $W^{m,p}(V)$ and
\[
D^\alpha\widetilde u=\widetilde{D^\alpha u},\qquad |\alpha|\leq m,
\qquad
\|\widetilde u\|_{W^{m,p}(V)}=\|u\|_{W^{m,p}(U)}.
\]
Moreover, $u\in W_0^{m,p}(U)$.
\end{lemma}
\begin{proof}
If $\varphi\in C_c^\infty(U)$, its extension by zero to $V$ is a smooth test function. For this test function, the identity defining $D^\alpha w$ on $V$ reduces to the identity on $U$. This identifies the derivatives of the restriction. The norm inequality follows by restricting each integral to $U$; no regularity of its boundary is required.

For extension by zero, the case $K=\varnothing$ is immediate. Choose $\eta\in C_c^\infty(U)$ equal to one on a neighborhood of $K$. Locality of weak derivatives gives $D^\alpha u=0$ almost everywhere on $U\setminus K$: it suffices to test the weak identity against functions in $C_c^\infty(U\setminus K)$ and use uniqueness of the function representing the zero distribution. If $\varphi\in C_c^\infty(V)$, then $\eta\varphi\restriction_U$ is an admissible test function and
\[
\begin{aligned}
\int_V\widetilde u\,D^\alpha\varphi\,dx
&=\int_Uu\,D^\alpha(\eta\varphi)\,dx\\
&=(-1)^{|\alpha|}\int_UD^\alpha u\,\eta\varphi\,dx\\
&=(-1)^{|\alpha|}\int_V\widetilde{D^\alpha u}\,\varphi\,dx.
\end{aligned}
\]
The first equality uses $\eta=1$ near $K$ and $u=0$ outside $K$; the last uses the same property for $D^\alpha u$. Each extension preserves the $L^p(V)$ norm, and therefore also the Sobolev norm.

To prove the final assertion, apply the preceding argument with $V=\mathbb R^n$. If $\rho_\varepsilon$ is a mollifier with support contained in $\overline B_{\mathrm{euc}}(0,\varepsilon)$, Fubini's theorem and the weak identity give
\[
D^\alpha(\rho_\varepsilon*\widetilde u)
=\rho_\varepsilon*\widetilde{D^\alpha u}.
\]
Approximation by mollifiers in $L^p(\mathbb R^n)$ implies convergence of each derivative in $L^p(\mathbb R^n)$. For sufficiently small $\varepsilon$, the support of the convolution is contained in $K+\overline B_{\mathrm{euc}}(0,\varepsilon)\Subset U$.
Their restrictions belong to $C_c^\infty(U)$ and converge to $u$ in $W^{m,p}(U)$, proving membership in $W_0^{m,p}(U)$.
\end{proof}

\begin{proposition}[Density of smooth functions]
\label{prop:densidad-sobolev-orden-entero}
Let $m\in\mathbb N_0$ and $1\leq p<\infty$. For every open subset $\Omega\subseteq\mathbb R^n$, the space $C^\infty(\Omega)\cap W^{m,p}(\Omega)$ is dense in $W^{m,p}(\Omega)$. Moreover,
\[
W_0^{m,p}(\mathbb R^n)=W^{m,p}(\mathbb R^n).
\]
Here $W^{0,p}(\Omega):=L^p(\Omega)$ and $W_0^{0,p}(\Omega)=L^p(\Omega)$, by Proposition~\ref{densidad de Cc en Lp}.
\end{proposition}
\begin{proof}
Successive applications of the weak Leibniz rule in Proposition~\ref{linealidad derivada debil} give, for $|\alpha|\leq m$,
\[
D^\alpha(\chi u)
=\sum_{\beta\leq\alpha}\binom{\alpha}{\beta}
(D^\beta\chi)D^{\alpha-\beta}u,
\]
where $\alpha,\beta\in\mathbb N_0^n$ and $\beta\leq\alpha$ means $\beta_i\leq\alpha_i$ for every $i\in\{1,\dots,n\}$. The induction producing this formula uses the first-order rule on each summand and Pascal's identity to combine the coefficients.

On $\mathbb R^n$, take the cutoffs $\chi_R$ from Proposition~\ref{Cc es denso en W1p}. The summand with $\beta=0$ converges to $D^\alpha u$ in $L^p(\mathbb R^n)$ by dominated convergence. If $\beta\neq0$, then
\[
\|(D^\beta\chi_R)D^{\alpha-\beta}u\|_{L^p(\mathbb R^n)}
\leq R^{-|\beta|}\|D^\beta\chi\|_{L^\infty(\mathbb R^n)}
\|D^{\alpha-\beta}u\|_{L^p(\mathbb R^n)}\longrightarrow0.
\]
Thus, $\chi_Ru\to u$ in $W^{m,p}(\mathbb R^n)$. For fixed $R$, Fubini's theorem and the definition of weak derivative give
\[
D^\alpha(\rho_\varepsilon*(\chi_Ru))
=\rho_\varepsilon*D^\alpha(\chi_Ru).
\]
Continuity of translations and the mollification estimate used in Proposition~\ref{Cc es denso en W1p} prove convergence in $L^p(\mathbb R^n)$ for each of the finitely many multi-indices. These mollifications are smooth and compactly supported, proving the stated equality.

For a general open subset, choose an exhaustion $\Omega_j\Subset\Omega_{j+1}\Subset\Omega$ with union $\Omega$, and a smooth partition of unity $(\eta_j)_{j\in\mathbb N}$ whose supports are compact and contained in the open sets $V_j:=\Omega_{j+1}\setminus\overline{\Omega_{j-1}}$; we may take $\Omega_0=\Omega_{-1}=\varnothing$ and adjust the first two members. The existence of this partition follows from Theorem~\ref{particionesdelaunidad} by grouping the functions subordinate to each $V_j$. The family $(V_j)_{j\in\mathbb N}$ is locally finite. The Leibniz formula shows that $\eta_ju\in W^{m,p}(V_j)$, and its support is contained in the compact set $\operatorname{supp}\eta_j\Subset V_j$. Lemma~\ref{lem:restriccion-extension-local-euclidiana} allows its extension by zero to $\mathbb R^n$ without changing either its norm or its weak derivatives.

Given $\delta>0$, for each $j$ choose $\varepsilon_j>0$ sufficiently small that
\[
v_j:=\rho_{\varepsilon_j}*(\eta_ju)\in C_c^\infty(V_j),
\qquad
\|v_j-\eta_ju\|_{W^{m,p}(\mathbb R^n)}<\delta2^{-j}.
\]
The sum $v:=\displaystyle\sum_{j=1}^{\infty}v_j$ is locally finite and therefore smooth on $\Omega$. The series of errors $\displaystyle\sum_{j=1}^{\infty}(v_j-\eta_ju)$ converges absolutely in the Banach space $W^{m,p}(\mathbb R^n)$. Its limit, restricted to $\Omega$, agrees with $v-u$: on each $\Omega_N$, only finitely many terms occur, by local finiteness of the supports. In particular, $v-u\in W^{m,p}(\Omega)$ before we take its norm, and the triangle inequality gives
\[
\|v-u\|_{W^{m,p}(\Omega)}
\leq\sum_{j=1}^{\infty}\|v_j-\eta_ju\|_{W^{m,p}(\mathbb R^n)}
\leq\delta.
\]
This proves that $v\in W^{m,p}(\Omega)$ and establishes density.
\end{proof}

\begin{proposition}[Extension by zero and derivatives of $W_0^{m,p}(\Omega)$]
\label{prop:extension-cero-W0-orden-entero}
Let $\Omega\subseteq\mathbb R^n$ be open, and let $m\in\mathbb N_0$ and $1\leq p<\infty$. Extension by zero $Z_\Omega$ is an isometry
\[
Z_\Omega\colon W_0^{m,p}(\Omega)\longrightarrow W^{m,p}(\mathbb R^n)
\]
and, for every multi-index $\alpha\in\mathbb N_0^n$ with $|\alpha|\leq m$,
\[
D^\alpha(Z_\Omega u)=Z_\Omega(D^\alpha u),
\qquad
D^\alpha u\in W_0^{m-|\alpha|,p}(\Omega).
\]
If $v\in W^{m,p}(\mathbb R^n)$ has support contained in a compact set $K\subseteq\Omega$, then $v\restriction_\Omega\in W_0^{m,p}(\Omega)$.
\end{proposition}
\begin{proof}
Choose $u_j\in C_c^\infty(\Omega)$ with $u_j\to u$ in $W^{m,p}(\Omega)$. For these functions, $D^\alpha(Z_\Omega u_j)=Z_\Omega(D^\alpha u_j)$ and
\[
\|Z_\Omega u_j-Z_\Omega u_\ell\|_{W^{m,p}(\mathbb R^n)}
=\|u_j-u_\ell\|_{W^{m,p}(\Omega)}.
\]
Completeness and Lemma~\ref{lema 16.13 generalizado} identify the limit with $Z_\Omega u$ and all its derivatives with the corresponding extensions by zero. Equality of norms passes to the limit. Moreover, $D^\alpha u_j\to D^\alpha u$ in $W^{m-|\alpha|,p}(\Omega)$; since each $D^\alpha u_j\in C_c^\infty(\Omega)$, membership in the indicated closed subspace follows.

For the final assertion, the mollifications of $v$ have support in $K+\overline B_{\mathrm{euc}}(0,\varepsilon)\Subset\Omega$ if $\varepsilon$ is sufficiently small, and converge in $W^{m,p}(\mathbb R^n)$ by the proof of Proposition~\ref{prop:densidad-sobolev-orden-entero}. Their restrictions approximate $v\restriction_\Omega$ by functions in $C_c^\infty(\Omega)$.
\end{proof}

The definition by closure requires no regularity of $\partial\Omega$. In particular, the preceding proposition does not identify $W_0^{m,p}(\Omega)$ with the kernel of a trace on an arbitrary open subset. For $W^{m,p}(\Omega)$, one must distinguish smooth approximation in the interior from extension across the boundary.

\begin{definition}[Sobolev extension domain]
\label{def:dominio-extension-sobolev}
We say that an open subset $\Omega$ is an extension domain for $W^{m,p}(\Omega)$ if there is a continuous linear operator
$\mathcal E\colon W^{m,p}(\Omega)\longrightarrow W^{m,p}(\mathbb R^n)$
such that $(\mathcal Eu)\restriction_\Omega=u$ almost everywhere.
\end{definition}
If $p<\infty$ and such an operator exists, approximation of $\mathcal Eu$ by functions in $C_c^\infty(\mathbb R^n)$ shows that their restrictions are dense in $W^{m,p}(\Omega)$. These restrictions need not vanish on the boundary; their density does not assert that $W_0^{m,p}(\Omega)=W^{m,p}(\Omega)$. The restriction to finite exponents in the equality on $\mathbb R^n$ is necessary: the constant function one belongs to $W^{m,\infty}(\mathbb R^n)$ and cannot be approximated uniformly by compactly supported functions.

\begin{definition}\label{def:sobolev-en-abiertos-euclidianos-define-esimo-exponente-critico-sobolev-note}\index{critical Sobolev exponent!higher-order}
If $m<n$ and $p\in [1,\displaystyle\frac{n}{m})$, the $m$th critical Sobolev exponent is defined by \[p_{m}^{*}:=\frac{np}{n-mp}.\] Note that $p_{1}^{*}=p^{*}$.
\end{definition}

The following formulation gathers the higher-order embeddings for finite exponents.

\begin{theorem}[Sobolev embedding with finite orders and exponents]
\label{prop:encaje-continuo-indices-finitos}
\label{sobolev generalizado}
\index{Sobolev embedding theorem!higher-order}
Let $\Omega\subseteq\mathbb R^n$ be open and bounded. Let $m>k\geq0$ be integers and $1\leq p,q<\infty$ satisfy
\[
m-\frac{n}{p}\geq k-\frac{n}{q}.
\]
Then $W_0^{m,p}(\Omega)\hookrightarrow W_0^{k,q}(\Omega)$ continuously. If $\Omega$ is an extension domain for $W^{m,p}(\Omega)$, then also $W^{m,p}(\Omega)\hookrightarrow W^{k,q}(\Omega)$ continuously.
\end{theorem}
\begin{proof}
For $q\leq p$, Hölder's inequality gives, for $|\alpha|\leq k$,
\[
\|D^\alpha u\|_{L^q(\Omega)}
\leq\lambda_n(\Omega)^{\frac{1}{q}-\frac{1}{p}}
\|D^\alpha u\|_{L^p(\Omega)}.
\]
Summing over the finitely many multi-indices gives $\|u\|_{W^{k,q}(\Omega)}\leq C\|u\|_{W^{m,p}(\Omega)}$.
If $u\in W_0^{m,p}(\Omega)$, applying this inequality to approximations of $C_c^\infty(\Omega)$ also shows that $u\in W_0^{k,q}(\Omega)$.
Assume $q>p$ and set $d:=m-k$ and
\[
\frac{1}{r_i}:=\frac{1}{p}-\frac{i}{d}\left(\frac{1}{p}-\frac{1}{q}\right),
\qquad i\in\{0,\dots,d\}.
\]
The exponents are finite, $r_0=p$, $r_d=q$, and
\[
\frac{1}{r_{i+1}}\geq\frac{1}{r_i}-\frac{1}{n}
\qquad (0\leq i<d).
\]
If $r_i<n$, this means $r_{i+1}\leq r_i^*$; if $r_i=n$ or $r_i>n$, the next exponent is still finite. In all three cases, the corresponding part of Theorem~\ref{desigualdad de poincaré} gives $W_0^{1,r_i}(\Omega)\hookrightarrow L^{r_{i+1}}(\Omega)$. Applying it to $D^\alpha\varphi$, with $|\alpha|\leq m-i-1$ and $\varphi\in C_c^\infty(\Omega)$, yields
\[
\|D^\alpha\varphi\|_{L^{r_{i+1}}(\Omega)}
\leq C_i\|D^\alpha\varphi\|_{W^{1,r_i}(\Omega)}.
\]
Summing over the multi-indices and using the equivalence of norms in finite dimensions, Proposition~\ref{normas en R^{n}}, gives
\[
\begin{aligned}
\|\varphi\|_{W^{m-i-1,r_{i+1}}(\Omega)}
&\leq C_i\left(
\sum_{|\alpha|\leq m-i-1}
\|D^\alpha\varphi\|_{W^{1,r_i}(\Omega)}^{r_{i+1}}
\right)^{\frac{1}{r_{i+1}}}\\
&\leq C_i'\|\varphi\|_{W^{m-i,r_i}(\Omega)}.
\end{aligned}
\]
In the last inequality, each derivative of order at most $m-i$ appears a finite number of times, depending only on $n,m,i$.

If $u\in W_0^{m-i,r_i}(\Omega)$, choose $\varphi_\nu\in C_c^\infty(\Omega)$ converging to $u$ in $W^{m-i,r_i}(\Omega)$. The preceding estimate, applied to $\varphi_\nu-\varphi_\mu$, shows that this sequence is Cauchy in $W^{m-i-1,r_{i+1}}(\Omega)$. By completeness, it converges to an element $v$ of this space, and $v\in W_0^{m-i-1,r_{i+1}}(\Omega)$ by the definition of closure. Since $\Omega$ has finite measure, convergence in $L^{r_i}(\Omega)$ and in $L^{r_{i+1}}(\Omega)$ implies convergence in $L^1(\Omega)$; uniqueness of this limit gives $v=u$ almost everywhere. Passing to the limit in the estimate yields a continuous embedding
\[
W_0^{m-i,r_i}(\Omega)\hookrightarrow
W_0^{m-i-1,r_{i+1}}(\Omega).
\]
The composition of these $d$ embeddings starts at $W_0^{m,p}(\Omega)$ and ends at $W_0^{m-d,r_d}(\Omega)=W_0^{k,q}(\Omega)$, as required.

If there is an extension operator $\mathcal E$, choose a ball $B$ containing $\overline\Omega$ and $\chi\in C_c^\infty(B)$ equal to one on a neighborhood of $\overline\Omega$. Set $Tu:=(\chi\mathcal Eu)\restriction_B$. The Leibniz rule and Proposition~\ref{prop:extension-cero-W0-orden-entero} give $Tu\in W_0^{m,p}(B)$, $Tu\restriction_\Omega=u$, and
\[
\begin{aligned}
\|Tu\|_{W^{m,p}(B)}
&\leq\|\chi\mathcal Eu\|_{W^{m,p}(\mathbb R^n)}\\
&\leq C_\chi\|\mathcal Eu\|_{W^{m,p}(\mathbb R^n)}
\leq C_\chi\|\mathcal E\|\,\|u\|_{W^{m,p}(\Omega)}.
\end{aligned}
\]
The assertion already proved on $B$ and continuity of restriction from Lemma~\ref{lem:restriccion-extension-local-euclidiana} yield
\[
\|u\|_{W^{k,q}(\Omega)}
\leq\|Tu\|_{W^{k,q}(B)}
\leq C_BC_\chi\|\mathcal E\|\,\|u\|_{W^{m,p}(\Omega)}.
\]
This proves the second embedding. We have not required $u$ to vanish on $\partial\Omega$: the cutoff vanishes near the boundary of the larger ball $B$.
\end{proof}

When $(m-k)p<n$, the Sobolev index is expressed through the critical exponent $p_{m-k}^*$. The usual formulation is a direct consequence of the preceding theorem.

\begin{corollary}[Sobolev embedding and critical exponent]
\label{cor:sobolev-exponente-critico-euclidiano}
\label{lema para sobolev generalizado}
Let $\Omega\subseteq\mathbb R^n$ be a bounded open subset. Let $j\geq0$ and $d\geq1$ be integers, and let $1\leq p<\displaystyle\frac{n}{d}$. Then
\[
W_0^{j+d,p}(\Omega)\hookrightarrow W_0^{j,q}(\Omega)
\qquad\text{for every }1\leq q\leq p_d^*
:=\frac{np}{n-dp},
\]
and these inclusions are continuous. In particular, $W_0^{d,p}(\Omega)\hookrightarrow L^q(\Omega)$ for the same exponents. If $\Omega$ is an extension domain for $W^{j+d,p}(\Omega)$, then also $W^{j+d,p}(\Omega)\hookrightarrow W^{j,q}(\Omega)$ continuously.
\end{corollary}
\begin{proof}
The condition $q\leq p_d^*$ is equivalent to
\[
j+d-\frac{n}{p}\geq j-\frac{n}{q}.
\]
Apply Theorem~\ref{sobolev generalizado} with total orders $m=j+d$ and $k=j$. For $j=0$, use $W_0^{0,q}(\Omega)=L^q(\Omega)$.
\end{proof}

\section{The Rellich--Kondrashov theorem in Euclidean spaces}
The Poincaré inequality provides continuous embeddings of $W_{0}^{1,p}(\Omega)$ into $L^{q}(\Omega)$ up to the critical exponent. Continuity, however, does not ensure that a bounded sequence has a convergent subsequence in the target space. The Rellich--Kondrashov theorem determines the exponents for which this gain in compactness occurs, a decisive feature in many existence arguments. The proof begins with an estimate for translations.
\begin{lemma}[Translation estimate in $L^p(\mathbb R^n)$]\label{lema 17.11}
Let $1\leq p<\infty$, $u\in W^{1,p}(\mathbb R^n)$, and $\xi\in\mathbb R^n$. If $T_\xi u(x):=u(x-\xi)$, then
\[
\|T_\xi u-u\|_{L^p(\mathbb R^n)}
\leq \|\xi\|\sum_{i=1}^{n}\|D_i u\|_{L^p(\mathbb R^n)}
\leq n^{1-\frac{1}{p}}\|\xi\|\|\nabla u\|_{(L^p(\mathbb R^n))^n}.
\]
\end{lemma}
\begin{proof}
For $\varphi\in C_c^\infty(\mathbb R^n)$, the fundamental theorem of calculus along the segment $t\mapsto x-t\xi$ gives
\[
\varphi(x-\xi)-\varphi(x)
=-\int_0^1\nabla\varphi(x-t\xi)\cdot\xi\,dt.
\]
Minkowski's integral inequality from Theorem~\ref{teo:minkowski-integral} and translation invariance of the measure imply
\[
\begin{aligned}
\|T_\xi\varphi-\varphi\|_{L^p(\mathbb R^n)}
&\leq\|\xi\|\int_0^1\sum_{i=1}^{n}
\|D_i\varphi(\,\cdot-t\xi)\|_{L^p(\mathbb R^n)}\,dt\\
&=\|\xi\|\sum_{i=1}^{n}\|D_i\varphi\|_{L^p(\mathbb R^n)}.
\end{aligned}
\]
For $p=1$, this also follows directly by applying Tonelli's theorem to the pointwise estimate. For $u\in W^{1,p}(\mathbb R^n)$, choose $\varphi_j\to u$ in $W^{1,p}(\mathbb R^n)$ by Proposition~\ref{Cc es denso en W1p}. Since translations are isometries,
\[
\|T_\xi u-u\|_{L^p(\mathbb R^n)}
\leq2\|u-\varphi_j\|_{L^p(\mathbb R^n)}+\|T_\xi\varphi_j-\varphi_j\|_{L^p(\mathbb R^n)}.
\]
Let $j\to\infty$. The second inequality in the statement is Hölder's inequality for the finite sum, with equality of the factors when $p=1$.
\end{proof}

\begin{theorem}[Rellich--Kondrashov with finite orders and exponents]
\label{rellich kondrashov generalizado}
\index{Rellich Kondrashov theorem@Rellich--Kondrashov theorem!higher-order}
Let $\Omega\subseteq\mathbb R^n$ be a bounded open subset. Let $m>k\geq0$ be integers and $1\leq p,q<\infty$ satisfy
\begin{equation}
\label{eq:rellich-indice-estricto-euclidiano}
m-\frac{n}{p}>k-\frac{n}{q}.
\end{equation}
Then the embedding
\[
W_0^{m,p}(\Omega)\hookrightarrow W_0^{k,q}(\Omega)
\]
is compact. If, in addition, $\Omega$ is an extension domain for $W^{m,p}(\Omega)$, the embedding
\[
W^{m,p}(\Omega)\hookrightarrow W^{k,q}(\Omega)
\]
is also compact. The continuity constants depend on $n,m,k,p,q,\Omega$ and, in the second assertion, on the norm of the chosen extension operator.
\end{theorem}
\begin{proof}
Let $(u_j)_{j\in\mathbb N}$ be a bounded sequence in $W_0^{m,p}(\Omega)$. We begin by proving convergence of a subsequence in $W^{m-1,p}(\Omega)$. Extend by zero and write $v_j:=Z_\Omega u_j$. By Proposition~\ref{prop:extension-cero-W0-orden-entero}, $(v_j)$ is bounded in $W^{m,p}(\mathbb R^n)$ and, for each $|\alpha|\leq m-1$,
\[
D^\alpha v_j\in W^{1,p}(\mathbb R^n),
\qquad
\operatorname{supp}(D^\alpha v_j)\subseteq\overline\Omega.
\]
If $M_0:=\displaystyle\sup_{j\in\mathbb N}
\|u_j\|_{W^{m,p}(\Omega)}<\infty$, Lemma~\ref{lema 17.11} gives
\[
\sup_{j\in\mathbb N}
\|T_\xi D^\alpha v_j-D^\alpha v_j\|_{L^p(\mathbb R^n)}
\leq C_{n,p}M_0\|\xi\|,
\qquad |\alpha|\leq m-1.
\]
Apply Proposition~\ref{corolario 14.47 relativamente compacto en Lp} on the open subset $\mathbb R^n$. The translation condition holds on every relatively compact open subset by the global estimate; the tail condition holds with zero error upon choosing a ball containing $\overline\Omega$. Boundedness in $L^p(\mathbb R^n)$ is already known. Therefore, $(D^\alpha v_j)_j$ has a subsequence converging in $L^p(\mathbb R^n)$.

Enumerate the finitely many multi-indices of order at most $m-1$, extract a subsequence for the first, then a subsequence of that one for the second, and so on. The last subsequence retains all the convergences. Restricting to $\Omega$ and applying Lemma~\ref{lema 16.13 generalizado} yields convergence in $W^{m-1,p}(\Omega)$; when $m=1$, completeness of $L^p(\Omega)$ suffices. Each term belongs to $W_0^{m-1,p}(\Omega)$ by its approximation in $W^{m,p}(\Omega)$ by test functions. This subspace is closed, so the limit belongs to it. The argument uses strong compactness in $L^p(\Omega)$ and also works for $p=1$.

If $q\leq p$, the continuous inclusion $W^{m-1,p}(\Omega)\hookrightarrow W^{k,q}(\Omega)$, obtained from the finite measure of $\Omega$ and $k\leq m-1$, proves the desired convergence. If $q>p$, the strict inequality \eqref{eq:rellich-indice-estricto-euclidiano} allows us to choose a finite exponent $r>q$ such that
\[
m-\frac{n}{p}\geq k-\frac{n}{r}.
\]
Theorem~\ref{prop:encaje-continuo-indices-finitos} shows that $(u_j)$ is bounded in $W^{k,r}(\Omega)$. For each $|\alpha|\leq k$, the interpolation inequality from Theorem~\ref{Ejercicio 14.67} gives
\[
\|D^\alpha(u_j-u_\ell)\|_{L^q(\Omega)}
\leq
\|D^\alpha(u_j-u_\ell)\|_{L^p(\Omega)}^\theta
\|D^\alpha(u_j-u_\ell)\|_{L^r(\Omega)}^{1-\theta},
\quad
\frac{1}{q}=\frac{\theta}{p}+\frac{1-\theta}{r},
\]
with $0<\theta<1$. Along the preceding subsequence, the first factor tends to zero and the second remains bounded. Summing over the multi-indices, we obtain a Cauchy sequence in $W^{k,q}(\Omega)$, which converges by completeness. In both cases, Theorem~\ref{prop:encaje-continuo-indices-finitos} places each term in $W_0^{k,q}(\Omega)$; closedness ensures that the limit also belongs to this space.

For $W^{m,p}(\Omega)$, use the operator $Tu=(\chi\mathcal Eu)\restriction_B\in W_0^{m,p}(B)$ constructed in the proof of Theorem~\ref{prop:encaje-continuo-indices-finitos}. For a bounded sequence $(u_j)_{j\in\mathbb N}$, with the same operator and the same cutoff, we have
\[
\sup_{j\in\mathbb N}\|Tu_j\|_{W^{m,p}(B)}
\leq C_\chi\|\mathcal E\|
\sup_{j\in\mathbb N}\|u_j\|_{W^{m,p}(\Omega)}<\infty.
\]
The first part provides a subsequence $Tu_{j_\nu}\to z$ in $W_0^{k,q}(B)$. Since $Tu_{j_\nu}\restriction_\Omega=u_{j_\nu}$,
\[
\|u_{j_\nu}-z\restriction_\Omega\|_{W^{k,q}(\Omega)}
\leq\|Tu_{j_\nu}-z\|_{W^{k,q}(B)}\longrightarrow0.
\]
Compactness is applied on the fixed ball $B$, not on all of $\mathbb R^n$.
\end{proof}

\begin{corollary}[Loss of one derivative with the same exponent]
\label{cor:rellich-descenso-un-orden}
For every bounded open subset $\Omega\subseteq\mathbb R^n$, every integer $m\geq1$, and every $1\leq p<\infty$, the embedding
\[
W_0^{m,p}(\Omega)\hookrightarrow W_0^{m-1,p}(\Omega)
\]
is compact. If $\Omega$ is an extension domain for $W^{m,p}(\Omega)$, then $W^{m,p}(\Omega)\hookrightarrow W^{m-1,p}(\Omega)$ is also compact.
\end{corollary}
\begin{proof}
Take $k=m-1$ and $q=p$ in the preceding theorem. The difference between the Sobolev indices is exactly one, in every dimension.
\end{proof}

The subcritical embeddings follow by choosing the corresponding exponents in the general theorem.

\begin{corollary}[First-order Rellich--Kondrashov]\label{rellich kondrashov 1}\index{Rellich Kondrashov theorem@Rellich--Kondrashov theorem!first-order}
If $\Omega\subseteq\mathbb{R}^{n}$ is open and bounded, $p\in[1,n)$, and $q\in [1,p^{*})$, then the inclusion $W_{0}^{1,p}(\Omega)\subseteq L^{q}(\Omega)$ is compact.
\end{corollary}
\begin{proof}
Take $m=1$ and $k=0$ in Theorem~\ref{rellich kondrashov generalizado}. Since $p<n$, the condition $1-\displaystyle\frac{n}{p}>-\displaystyle\frac{n}{q}$ is equivalent to $q<p^*$.
\end{proof}
At higher orders, the finite critical exponent appears when the product of the increase in order and $p$ is less than $n$.
\begin{corollary}[Rellich--Kondrashov and critical exponent]
\label{cor:rellich-exponente-critico-orden-superior}
\index{Rellich Kondrashov theorem@Rellich--Kondrashov theorem!higher-order}
Let $\Omega\subseteq\mathbb R^n$ be a bounded open subset, let $j\in\mathbb N_0$ and $m\in\{1,\dots,n-1\}$, and let $1\leq p<\displaystyle\frac{n}{m}$. Then, for every $1\leq q<p_m^*$, the inclusion
\[
W_0^{j+m,p}(\Omega)\hookrightarrow W_0^{j,q}(\Omega)
\]
is compact. At the endpoint $q=p_m^*$, the embedding is continuous by Theorem~\ref{sobolev generalizado}.
\end{corollary}

\begin{proof}
Since $mp<n$, we have
\[
q<p_m^*
\quad\Longleftrightarrow\quad
j+m-\frac{n}{p}>j-\frac{n}{q}.
\]
Apply Theorem~\ref{rellich kondrashov generalizado} with total orders $j+m$ and $j$. Continuity at the endpoint is Corollary~\ref{cor:sobolev-exponente-critico-euclidiano}.
\end{proof}

\begin{remark}[The critical endpoint and assumptions on the domain]
\label{obs:limites-rellich-euclidiano}
The strict inequality cannot in general be replaced by equality. Let $m>k\geq0$, $(m-k)p<n$, $q=\displaystyle\frac{np}{n-(m-k)p}$, and $B_{\mathrm{euc}}(x_0,r)\Subset\Omega$. Choose a nonzero $\varphi\in C_c^\infty(B_{\mathrm{euc}}(0,1))$ and set, for $0<\varepsilon<r$,
\[
u_\varepsilon(x):=\varepsilon^{m-\frac{n}{p}}
\varphi\left(\frac{x-x_0}{\varepsilon}\right).
\]
For each $|\alpha|\leq m$,
\[
\|D^\alpha u_\varepsilon\|_{L^p(\Omega)}
=\varepsilon^{m-|\alpha|}\|D^\alpha\varphi\|_{L^p(B_{\mathrm{euc}}(0,1))}.
\]
Thus, the family is bounded in $W_0^{m,p}(\Omega)$. For some $|\beta|=k$, the function $D^\beta\varphi$ is nonzero and
\[
\|D^\beta u_\varepsilon\|_{L^q(\Omega)}
=\|D^\beta\varphi\|_{L^q(B_{\mathrm{euc}}(0,1))}>0.
\]
However, $D^\beta u_\varepsilon\to0$ almost everywhere as $\varepsilon\to 0^{+}$. No subsequence can converge in $L^q(\Omega)$. This prevents compactness at the critical exponent, although continuity does hold.

Nor does boundedness of an open subset suffice for compactness on $W^{m,p}(\Omega)$ without the zero subscript. Let $\Omega$ be the union of disjoint open balls $(B_j)_{j\in\mathbb N}$ contained in a fixed ball and accumulating only at a point outside $\Omega$. The functions
\[
u_j:=\lambda_n(B_j)^{-\frac{1}{p}}\mathbf1_{B_j}
\]
are constant on each component of $\Omega$, so all their positive-order weak derivatives vanish and $\|u_j\|_{W^{m,p}(\Omega)}=1$. For $j\neq\ell$, we have $\|u_j-u_\ell\|_{L^p(\Omega)}=2^{\displaystyle\frac{1}{p}}$. Thus, even the loss of one derivative at the same exponent does not yield compactness on this $W^{m,p}(\Omega)$.

The condition $m>k$ is also essential: on a ball of finite measure, the index inequality holds with $m=k=0$ and $q<p$, but the inclusion $L^p(\Omega)\hookrightarrow L^q(\Omega)$ is not compact. On a cube, functions alternating between $1$ and $-1$ on dyadic intervals of the first coordinate are bounded in $L^p(\Omega)$, and any two differ on a set whose measure is half that of the cube. Their distances in $L^q(\Omega)$ are therefore bounded away from zero.
\end{remark}

\subsection{Global embeddings and localization in $\mathbb R^n$}

\begin{proposition}[Global continuous embeddings with finite exponents]
\label{prop:encajes-globales-Rn-finitos}
Let $m>k\geq0$ be integers and $1\leq p\leq q<\infty$ satisfy $m-\displaystyle\frac{n}{p}\geq k-\displaystyle\frac{n}{q}$. Then
\[
W^{m,p}(\mathbb R^n)\hookrightarrow W^{k,q}(\mathbb R^n)
\]
continuously. The same conclusion holds on an extension domain for $W^{m,p}(\Omega)$ by restricting its extension operator.
\end{proposition}
\begin{proof}
Let $Q_a:=a+(0,1)^n$, with $a\in\mathbb Z^n$, and $\chi_a(x):=\chi(x-a)$, where $\chi\in C_c^\infty((-1,2)^n)$ equals one on $[0,1]^n$. Proposition~\ref{prop:extension-cero-W0-orden-entero} and the Leibniz rule give $\chi_au\in W_0^{m,p}(a+(-1,2)^n)$. Applying Theorem~\ref{prop:encaje-continuo-indices-finitos} on these translated cubes yields
\[
\|u\|_{W^{k,q}(Q_a)}
\leq C\|u\|_{W^{m,p}(a+(-1,2)^n)},
\]
with the same constant for every $a$. The cubes $Q_a$ cover $\mathbb R^n$ except for a null set, and the enlarged cubes have multiplicity at most $3^n$. Since $p\leq q$, for every nonnegative family $(b_a)$ we have $\displaystyle(\displaystyle\sum_{a\in\mathbb Z^n} b_a^q)^{\displaystyle\frac{1}{q}}\leq(\displaystyle\sum_{a\in\mathbb Z^n} b_a^p)^{\displaystyle\frac{1}{p}}$. Consequently,
\[
\begin{aligned}
\|u\|_{W^{k,q}(\mathbb R^n)}
&\leq C\left(\sum_{a\in\mathbb Z^n}
\|u\|_{W^{m,p}(a+(-1,2)^n)}^p\right)^{\frac{1}{p}}\\
&\leq C3^{\frac{n}{p}}\|u\|_{W^{m,p}(\mathbb R^n)}.
\end{aligned}
\]
Composition with an extension operator followed by restriction proves the final assertion.
\end{proof}

\begin{corollary}[Compactness with common support and local compactness]
\label{cor:rellich-local-Rn}
Let $m>k\geq0$, $1\leq p,q<\infty$, and $m-\displaystyle\frac{n}{p}>k-\displaystyle\frac{n}{q}$. For each compact set $K\subseteq\mathbb R^n$, the inclusion into $W^{k,q}(\mathbb R^n)$ of the subspace
\[
W_K^{m,p}(\mathbb R^n)
:=\{u\in W^{m,p}(\mathbb R^n)\mid\operatorname{supp}u\subseteq K\}
\]
is compact. For any open subset $\Omega$, every bounded sequence in $W^{m,p}(\Omega)$ has a subsequence converging in $W^{k,q}(U)$ for each open subset $U\Subset\Omega$.
\end{corollary}
\begin{proof}
Let $(u_\nu)_{\nu\in\mathbb N}$ be a bounded sequence in
$W_K^{m,p}(\mathbb R^n)$, and choose a ball $B$ such that $K\Subset B$.
By Proposition~\ref{prop:extension-cero-W0-orden-entero},
$u_\nu|_B\in W_0^{m,p}(B)$ and
\[
 \|u_\nu|_B\|_{W^{m,p}(B)}
 =\|u_\nu\|_{W^{m,p}(\mathbb R^n)}.
\]
Theorem~\ref{rellich kondrashov generalizado} gives a subsequence
$(u_{\nu_a}|_B)$ and $v\in W_0^{k,q}(B)$ such that
$u_{\nu_a}|_B\to v$ in $W^{k,q}(B)$. If $Z_B$ denotes zero
extension, the same proposition applied at order $k$ gives
\[
 \|u_{\nu_a}-Z_Bv\|_{W^{k,q}(\mathbb R^n)}
 =\|u_{\nu_a}|_B-v\|_{W^{k,q}(B)}
 \longrightarrow0.
\]
The bound for the Rellich inclusion and this equality also prove
continuity of the inclusion under consideration. Thus it is a
compact operator.

Now let $(u_\nu)_{\nu\in\mathbb N}$ be bounded in $W^{m,p}(\Omega)$.
If $\Omega=\mathbb R^n$, take $\Omega_j=B(0,j)$. If
$\Omega\neq\mathbb R^n$, we may take
\[
 \Omega_j
 :=\left\{x\in\Omega:
       |x|<j,\ 
       \operatorname{dist}(x,\mathbb R^n\setminus\Omega)>
       \frac{1}{j}\right\},
 \qquad j\geq1.
\]
In either case $\overline{\Omega_j}\subseteq\Omega_{j+1}$,
$\overline{\Omega_j}$ is compact, and
$\Omega=\displaystyle\bigcup_{j\geq1}\Omega_j$.

Choose $\chi_j\in C_c^\infty(\Omega_{j+1})$ equal to one on a
neighborhood of $\overline{\Omega_j}$, and put
$K_j=\operatorname{supp}\chi_j\Subset\Omega$.
The Leibniz rule for weak derivatives gives, for $|\alpha|\leq m$,
\[
 D^\alpha(\chi_j u_\nu)
 =\displaystyle\sum_{\beta\leq\alpha}
       {\alpha\choose\beta}
       (D^\beta\chi_j)D^{\alpha-\beta}u_\nu.
\]
Each derivative of $\chi_j$ appearing here is bounded, and the
sum is finite. Lemma~\ref{lem:restriccion-extension-local-euclidiana}
allows us to extend the product by zero to $\mathbb R^n$.
Denoting that extension by $w_{\nu,j}$, we obtain
\[
 w_{\nu,j}\in W_{K_j}^{m,p}(\mathbb R^n),\qquad
 \|w_{\nu,j}\|_{W^{m,p}(\mathbb R^n)}
 \leq C_j\|u_\nu\|_{W^{m,p}(\Omega)}.
\]
The constant $C_j$ depends on the chosen cutoffs and on $m,p,n$,
but not on $\nu$.

The first part of the proof applies to $(w_{\nu,1})_\nu$.
There is an infinite set $I_1\subseteq\mathbb N$ such that
$w_{\nu,1}$, for $\nu\in I_1$ in increasing order, converges in
$W^{k,q}(\mathbb R^n)$. Applying it successively gives infinite sets
\[
 I_1\supseteq I_2\supseteq I_3\supseteq\cdots
\]
and functions $v_j\in W^{k,q}(\mathbb R^n)$ such that
$w_{\nu,j}\to v_j$ as $\nu\to\infty$ within $I_j$.
Choose recursively $\nu_a\in I_a$ with $\nu_a>\nu_{a-1}$.
For each fixed $j$, all terms with $a\geq j$ have index in $I_j$;
hence
\[
 w_{\nu_a,j}\longrightarrow v_j
 \quad\text{in }W^{k,q}(\mathbb R^n).
\]
Since $\chi_j=1$ on $\Omega_j$, continuity of restriction gives
$u_{\nu_a}|_{\Omega_j}\to v_j|_{\Omega_j}$ in
$W^{k,q}(\Omega_j)$.
If $j\leq\ell$, the same sequence restricted to $\Omega_j$ converges
both to $v_j|_{\Omega_j}$ and to $v_\ell|_{\Omega_j}$.
Uniqueness of the limit in $L^q(\Omega_j)$ implies
$v_j=v_\ell$ almost everywhere on $\Omega_j$.

These limits define a measurable function $u_\infty$ on $\Omega$:
with $\Omega_0=\varnothing$, take $u_\infty=v_j$ on
$\Omega_j\setminus\Omega_{j-1}$. The compatibility just proved
implies $u_\infty=v_j$ almost everywhere on $\Omega_j$.
Thus $u_\infty\in W^{k,q}_{\mathrm{loc}}(\Omega)$:
each compactly supported test function is contained in some
$\Omega_j$, where the weak derivative identities are those of $v_j$.
Finally, if $U\Subset\Omega$, its closure is contained in some
$\Omega_j$, and
\[
 \|u_{\nu_a}-u_\infty\|_{W^{k,q}(U)}
 \leq\|w_{\nu_a,j}-v_j\|_{W^{k,q}(\mathbb R^n)}
 \longrightarrow0.
\]
The same subsequence works for all such open sets $U$.
\end{proof}

\begin{remark}[Translations toward infinity]
\label{obs:no-compacidad-global-Rn}
The equality $W_0^{m,p}(\mathbb R^n)=W^{m,p}(\mathbb R^n)$ for $p<\infty$ does not make the global continuous embeddings compact. If $0\neq\varphi\in C_c^\infty(\mathbb R^n)$ and we choose points $x_j$ so that the supports of $u_j(x):=\varphi(x-x_j)$ are disjoint, the $W^{m,p}(\mathbb R^n)$ norms remain constant and
\[
\|u_j-u_\ell\|_{L^q(\mathbb R^n)}
=2^{\frac{1}{q}}\|\varphi\|_{L^q(\mathbb R^n)}
\qquad(j\neq\ell, 1\leq q<\infty).
\]
For $q=\infty$, the distance is $\|\varphi\|_{L^\infty(\mathbb R^n)}>0$. Global compactness fails even with a strict gain in indices. Localization or a common compact support prevents this escape.

The condition $q\geq p$ in the global continuous embeddings must also be retained. If $1\leq q<p<\infty$, choose $a$ with $\displaystyle\frac{n}{p}<a\leq \displaystyle\frac{n}{q}$. The function $u(x)=(1+|x|^2)^{-\displaystyle\frac{a}{2}}$ and all its derivatives of every order belong to $L^p(\mathbb R^n)$, whereas $u\notin L^q(\mathbb R^n)$. Inclusions based on finite measure apply to bounded domains or fixed supports, not to the whole Euclidean space.
\end{remark}

\section{Sobolev embeddings into Hölder spaces}

The Sobolev embeddings studied so far improve the integrability of a weakly differentiable function. In the supercritical regime, the gain is stronger: when the number of derivatives and the integrability exponent satisfy $mp>n$, functions in $W^{m,p}(\Omega)$ admit representatives with Hölder continuous classical derivatives. This property is fundamental in regularity theory. The embedding does not itself produce new derivatives: one must first obtain a priori estimates placing the weak solution in a sufficiently high space $W^{m,p}(\Omega)$; Morrey's embedding then converts this Sobolev regularity into quantitative classical regularity. The development in this section follows Morrey's strategy as presented in \cite[Section~12.3]{Leoni2017}.

\begin{definition}[Functions of class $C^k$ and $C^{k,\alpha}$ on the closure of an open subset]
\label{def:espacios-holder-euclidianos}
\index{Holder space@Hölder space}
Let $\Omega\subseteq\mathbb R^n$ be an open subset. We begin with the usual convention
\[
C^0(\overline\Omega)
:=
\left\{
 v\colon \overline\Omega\longrightarrow\mathbb R
 \middle|
 v\text{ is continuous}
\right\}.
\]
For $k\in\mathbb N\cup\{\infty\}$, write
\[
C^k(\overline\Omega)
:=
\left\{
 v\in C^0(\overline\Omega)
 \middle|
 \begin{array}{l}
 \text{there are an open subset $V\subseteq\mathbb R^n$, with
 $\overline\Omega\subseteq V$, and a function}\\[-1mm]
 \widetilde v\in C^k(V)\text{ such that }
 \widetilde v\restriction_{\overline\Omega}=v
 \end{array}
\right\}.
\]
Continuity shows that an extension initially defined on $\Omega$ determines, whenever it exists, a unique function on $\overline\Omega$. In applications, we shall therefore identify a function on $\Omega$ with its continuous representative on the closure. Moreover, for $k\in\mathbb N_0\cup\{\infty\}$, define
\[
C_c^k(\overline\Omega)
:=
\left\{
 v\in C^k(\overline\Omega)
 \middle|
 \operatorname{supp}_{\overline\Omega}(v)\text{ is compact}
\right\},
\]
where the support is taken in the topological space $\overline\Omega$. Thus, for $k\geq1$, the notation $C^k(\overline\Omega)$ does not merely mean that the derivatives have limits at the boundary. It requires the function to be the restriction of a function of class $C^k$ defined on an open neighborhood of $\overline\Omega$. This is also the convention used in trace theory.

Now let $K\subseteq\mathbb R^n$ be nonempty and let $w\colon K\longrightarrow\mathbb R$ be a bounded function. For $0<\alpha<1$, set
\[
[w]_{C^{0,\alpha}(K)}
:=
\sup_{\substack{x,y\in K\\x\neq y}}
\frac{|w(x)-w(y)|}{|x-y|^\alpha},
\qquad
\|w\|_{C^{0,\alpha}(K)}
:=
\|w\|_{C^0(K)}+[w]_{C^{0,\alpha}(K)},
\qquad
\|w\|_{C^0(K)}:=\sup_{x\in K}|w(x)|,
\]
with the convention that the supremum is zero when $K$ consists of a single point. For orders $k\geq1$, we also assume $K=\overline{\operatorname{Int}(K)}$. This condition ensures that derivatives on $K$ are determined by the restriction of the function: two extensions agreeing on $K$ have the same derivatives on its interior and, by continuity, on $K$. If $k\in\mathbb N$ and $w$ is of class $C^k$ on a neighborhood of $K$, define the intrinsic norm
\begin{equation}
\label{eq:norma-holder-intrinseca-sobre-conjunto}
\|w\|_{C^{k,\alpha}_{\mathrm{int}}(K)}
:=
\sum_{|\beta|\leq k}\|D^\beta w\|_{C^0(K)}
+
\sum_{|\beta|=k}[D^\beta w]_{C^{0,\alpha}(K)}.
\end{equation}
For $k=0$, write $\|w\|_{C^{0,\alpha}_{\mathrm{int}}(K)}:=
\|w\|_{C^{0,\alpha}(K)}$. Denote by $C^{k,\alpha}_{\mathrm{int}}(K)$ the space of restrictions to $K$ of functions of class $C^k$ on a neighborhood of $K$ for which this norm is finite. The subscript $\mathrm{int}$ recalls that the norm compares only the values of the function and its derivatives on $K$; it does not require an extension with global control.

Let $k\in\mathbb N_0$ and $0<\alpha<1$. For $v\in C^k(\mathbb R^n)$, define
\[
[v]_{C^{k,\alpha}(\mathbb R^n)}
:=
\sum_{|\beta|=k}[D^\beta v]_{C^{0,\alpha}(\mathbb R^n)}
\]
and
\[
\|v\|_{C^{k,\alpha}(\mathbb R^n)}
:=
\sum_{|\beta|\leq k}
\|D^\beta v\|_{L^\infty(\mathbb R^n)}
+
[v]_{C^{k,\alpha}(\mathbb R^n)}.
\]
Denote by $C^{k,\alpha}(\mathbb R^n)$ the space of functions for which this norm is finite.

Finally, suppose that $\Omega$ is bounded. For $v\in C^k(\overline\Omega)$, let
\[
\mathcal E_{\Omega}^{k,\alpha}(v)
:=
\left\{
 \widetilde v\in C^{k,\alpha}(\mathbb R^n)
 \middle|
 \widetilde v\restriction_{\overline\Omega}=v
\right\}.
\]
Define
\[
C^{k,\alpha}(\overline\Omega)
:=
\left\{
 v\in C^k(\overline\Omega)
 \middle|
 \mathcal E_{\Omega}^{k,\alpha}(v)\neq\varnothing
\right\}
\]
and equip this space with the extension norm
\begin{equation}
\label{eq:norma-holder-cerradura-por-extensiones}
\|v\|_{C^{k,\alpha}(\overline\Omega)}
:=
\inf_{\widetilde v\in\mathcal E_{\Omega}^{k,\alpha}(v)}
\|\widetilde v\|_{C^{k,\alpha}(\mathbb R^n)}.
\end{equation}
In particular, $C^{k,\alpha}(\overline\Omega)\subseteq C^k(\overline\Omega)$.
\end{definition}

\begin{remark}[Local and global formulations by extensions]
\label{obs:equivalencia-extensiones-holder-cerradura}
The preceding definition admits a local formulation. Indeed, a function $v$ belongs to $C^{k,\alpha}(\overline\Omega)$ if and only if there are open sets
\[
\overline\Omega\subseteq V_0\Subset V_1\Subset V
\]
and an extension $\widehat v\in C^k(V)$ of $v$, that is, $\widehat v\restriction_{\overline\Omega}=v$, such that
\[
\|\widehat v\|_{C^{k,\alpha}_{\mathrm{int}}(\overline{V_1})}<\infty.
\]
A global extension satisfies this condition immediately. Conversely, choose a cutoff function $\chi\in C_c^\infty(V_1)$ equal to one on a neighborhood of $\overline{V_0}$. The extension by zero of $\chi\widehat v$ outside $V_1$, denoted by $\widetilde{\chi\widehat v}$, belongs to $C^{k,\alpha}(\mathbb R^n)$ and agrees with $v$ on $\overline\Omega$. The intermediate open set $V_1$ ensures that all derivatives involved in the Leibniz rule are controlled on a compact subset of the domain of $\widehat v$. For $k=\infty$, the part of this formulation concerning the existence of a smooth extension recovers exactly the convention used for $C^\infty(\overline\Omega)$ in trace theory.

When $k=0$, the class thus defined agrees with the usual class of Hölder continuous functions on the compact set $\overline\Omega$. Indeed, if $H=[v]_{C^{0,\alpha}(\overline\Omega)}$, McShane's formula
\[
\widehat v(x)
:=
\inf_{y\in\overline\Omega}
\bigl(v(y)+H|x-y|^\alpha\bigr)
\]
extends $v$ to a Hölder continuous function of exponent $\alpha$ on $\mathbb R^n$, since $(x,y)\mapsto|x-y|^\alpha$ is a metric for $0<\alpha\leq1$. Multiplication by a cutoff function equal to one near $\overline\Omega$ then produces a bounded extension. Therefore, when $k=0$, the extension norm in \eqref{eq:norma-holder-cerradura-por-extensiones} is equivalent to the usual intrinsic norm on $\overline\Omega$.

For $k\geq1$, the extension condition is part of the convention. In the applications in this section, the functions come from $W_0^{m,p}(\Omega)$, and their extensions by zero provide the required extensions. In the geometric applications, the coordinate domains are balls and the localized functions vanish near the coordinate boundary.
\end{remark}

\begin{remark}[Changing the Hölder exponent]
\label{obs:cambio-exponente-holder-cerradura}
If $0<\alpha\leq\beta<1$ and $\Omega$ is bounded, then
\[
C^{k,\beta}(\overline\Omega)
\hookrightarrow
C^{k,\alpha}(\overline\Omega)
\]
continuously. To see this, fix a cutoff function $\chi\in C_c^\infty(\mathbb R^n)$ equal to one on a neighborhood of $\overline\Omega$. If $\widetilde v$ is an extension of $v$, then $\chi\widetilde v$ is still an extension, and the Leibniz rule gives
\[
\|\chi\widetilde v\|_{C^{k,\beta}(\mathbb R^n)}
\leq C_\chi\|\widetilde v\|_{C^{k,\beta}(\mathbb R^n)}.
\]
Since $\chi\widetilde v$ has support in a fixed compact set, each of its derivatives of order $k$ satisfies
\[
[w]_{C^{0,\alpha}(\mathbb R^n)}
\leq
C_{\chi,\alpha,\beta}
\bigl(\|w\|_{L^\infty(\mathbb R^n)}
+[w]_{C^{0,\beta}(\mathbb R^n)}\bigr).
\]
Taking the infimum over extensions $\widetilde v$ proves continuity of the embedding. We shall also use the fact that a function of class $C^1$ with bounded derivative is Lipschitz, and hence Hölder continuous with every exponent $0<\alpha<1$ on sets of bounded diameter.
\end{remark}

The quantity defined in \eqref{eq:norma-holder-cerradura-por-extensiones} is a norm. Homogeneity and the triangle inequality follow by taking extensions and adding them. Moreover,
\[
\|v\|_{C^0(\overline\Omega)}
\leq
\|v\|_{C^{k,\alpha}(\overline\Omega)},
\]
since the same inequality holds for every extension. The norm therefore vanishes only when $v=0$.

\begin{proposition}[Completeness of Hölder spaces]
\label{prop:completitud-holder-euclidiano}
For $k\in\mathbb N_0$ and $0<\alpha<1$, the space $C^{k,\alpha}(\mathbb R^n)$ is a Banach space. If $\Omega\subseteq\mathbb R^n$ is open and bounded, then $C^{k,\alpha}(\overline\Omega)$, with the norm \eqref{eq:norma-holder-cerradura-por-extensiones}, is also a Banach space.
\end{proposition}

\begin{proof}
Let $(v_j)_{j\in\mathbb N}$ be a Cauchy sequence in $C^{k,\alpha}(\mathbb R^n)$. For each multi-index $\beta$ with $|\beta|\leq k$, the sequence $(D^\beta v_j)$ converges uniformly to a bounded continuous function $v_\beta$. If $|\beta|<k$, the fundamental theorem of calculus along segments allows passage to the limit in
\[
D^\beta v_j(x+te_i)-D^\beta v_j(x)
=
\int_0^tD^{\beta+e_i}v_j(x+se_i)\,ds.
\]
It follows that $v_0\in C^k(\mathbb R^n)$ and $D^\beta v_0=v_\beta$ for $|\beta|\leq k$. If $|\beta|=k$, passing to the limit in the Hölder inequalities gives $v_\beta\in C^{0,\alpha}(\mathbb R^n)$. Finally, letting one index tend to infinity in the Cauchy condition yields $v_j\to v_0$ in $C^{k,\alpha}(\mathbb R^n)$.

Now let $(u_j)_{j\in\mathbb N}$ be a Cauchy sequence in $C^{k,\alpha}(\overline\Omega)$. Extract a subsequence $(u_{j_\nu})_{\nu\in\mathbb N}$ such that
\[
\|u_{j_{\nu+1}}-u_{j_\nu}\|_{C^{k,\alpha}(\overline\Omega)}
<2^{-\nu-1}.
\]
Choose $U_1\in\mathcal E_\Omega^{k,\alpha}(u_{j_1})$. For each $\nu\in\mathbb N$, the definition of the infimum gives $G_\nu\in\mathcal E_\Omega^{k,\alpha}
(u_{j_{\nu+1}}-u_{j_\nu})$ such that
\[
\|G_\nu\|_{C^{k,\alpha}(\mathbb R^n)}<2^{-\nu}.
\]
The series
\[
U_1+\sum_{\nu=1}^{\infty}G_\nu
\]
converges in $C^{k,\alpha}(\mathbb R^n)$ to a function $U$. Upon restriction to $\overline\Omega$, the partial sum of order $N$ is $u_{j_{N+1}}$. Therefore, $u:=U\restriction_{\overline\Omega}$ belongs to $C^{k,\alpha}(\overline\Omega)$, and $u_{j_\nu}\to u$ in the extension norm. Since the original sequence is Cauchy, the whole sequence converges to $u$.
\end{proof}

The central analytic ingredient is the following local estimate. Its advantage is that the power of $r$ is determined by the scaling of the problem.

\begin{lemma}[Morrey's estimate on cubes]
\label{lem:estimacion-morrey-cubos}
Let $n<p<\infty$ and let $Q\subseteq\mathbb R^n$ be a cube with sides parallel to the coordinate axes and side length $r>0$. There is a constant $C=C(n,p)>0$ such that
\begin{equation}
\label{eq:morrey-cubo}
|v(x)-v(y)|
\leq
C r^{1-\frac{n}{p}}
\|\nabla v\|_{L^p(Q)}
\end{equation}
for every $v\in C^1(\overline Q)$ and all $x,y\in Q$.
\end{lemma}

\begin{proof}
Denote by
\[
v_Q:=\frac{1}{\lambda_n(Q)}\int_Qv(z)\,dz
\]
the average of $v$ on $Q$. Fix $x\in Q$. Since $Q$ is convex, for every $y\in Q$ the segment joining $x$ to $y$ is contained in $Q$ and
\[
v(y)-v(x)
=
\int_0^1
\langle\nabla v(x+t(y-x)),y-x\rangle\,dt.
\]
Averaging with respect to $y$ and making, for each $t\in(0,1)$, the change of variables $z=x+t(y-x)$, we obtain
\[
|v(x)-v_Q|
\leq
\frac{1}{r^n}
\int_0^1
\int_{x+t(Q-x)}
|x-z|\,t^{-n-1}|\nabla v(z)|\,dz\,dt.
\]
If $z\in x+t(Q-x)$, then $|z-x|\leq\sqrt n\,rt$ and therefore $t\geq \displaystyle\frac{|z-x|}{\sqrt n\,r}$. Tonelli's theorem allows us to change the order of integration and gives
\[
|v(x)-v_Q|
\leq
C_n\int_Q\frac{|\nabla v(z)|}{|x-z|^{n-1}}\,dz.
\]
Let $p'=\displaystyle\frac{p}{p-1}$. Hölder's inequality and the inclusion $Q\subseteq B_{\mathrm{euc}}(x,\sqrt n\,r)$ imply
\[
|v(x)-v_Q|
\leq
C_n
\left(
\int_{B_{\mathrm{euc}}(x,\sqrt n\,r)}
|x-z|^{(1-n)p'}\,dz
\right)^{\frac{1}{p'}}
\|\nabla v\|_{L^p(Q)}.
\]
The radial integral is finite precisely because $p>n$, and its power of $r$ is $1-\displaystyle\frac{n}{p}$. Consequently,
\[
|v(x)-v_Q|
\leq
C r^{1-\frac{n}{p}}\|\nabla v\|_{L^p(Q)}.
\]
Applying the same inequality to $y$ and using $|v(x)-v(y)|\leq|v(x)-v_Q|+|v(y)-v_Q|$ gives \eqref{eq:morrey-cubo}.
\end{proof}

\begin{theorem}[First-order Morrey embedding]
\label{teo:morrey-primer-orden-W0}
\index{Morrey theorem@Morrey's theorem!first-order}
Let $\Omega\subseteq\mathbb R^n$ be a bounded open subset and let $n<p<\infty$. If
\[
\mu:=1-\frac{n}{p},
\]
then every element of $W_0^{1,p}(\Omega)$ has a unique representative in $C^{0,\mu}(\overline\Omega)$, and there is a constant $C=C(n,p,\Omega)>0$ such that
\begin{equation}
\label{eq:morrey-primer-orden-W0}
\|u\|_{C^{0,\mu}(\overline\Omega)}
\leq
C\|u\|_{W^{1,p}(\Omega)}.
\end{equation}
\end{theorem}

\begin{proof}
First let $\varphi\in C_c^\infty(\Omega)$, and denote by $\widetilde\varphi$ its extension by zero to $\mathbb R^n$. Given $x,y\in\mathbb R^n$, choose a cube of side length $2|x-y|$ containing both points. Lemma~\ref{lem:estimacion-morrey-cubos}, applied to $\widetilde\varphi$, gives
\[
|\widetilde\varphi(x)-\widetilde\varphi(y)|
\leq
C|x-y|^\mu
\|\nabla\widetilde\varphi\|_{L^p(\mathbb R^n)}.
\]
Therefore,
\[
[\widetilde\varphi]_{C^{0,\mu}(\mathbb R^n)}
\leq
C\|\nabla\varphi\|_{L^p(\Omega)}.
\]
Take a fixed cube $Q_0$ containing $\overline\Omega$ in its interior and a point $y_0\in Q_0\setminus\overline\Omega$. Since $\widetilde\varphi(y_0)=0$, the preceding estimate also gives
\[
\|\widetilde\varphi\|_{L^\infty(\mathbb R^n)}
\leq
C_{\Omega}\|\nabla\varphi\|_{L^p(\Omega)}.
\]
Thus, \eqref{eq:morrey-primer-orden-W0} holds for every $\varphi\in C_c^\infty(\Omega)$.

Now let $u\in W_0^{1,p}(\Omega)$ and choose $(\varphi_j)_{j\in\mathbb N}\subseteq C_c^\infty(\Omega)$ such that $\varphi_j\to u$ in $W^{1,p}(\Omega)$. The global inequality already proved shows that the sequence of extensions by zero $(\widetilde\varphi_j)_{j\in\mathbb N}$ is Cauchy in $C^{0,\mu}(\mathbb R^n)$. By Proposition~\ref{prop:completitud-holder-euclidiano}, it converges to a function $\widetilde u\in C^{0,\mu}(\mathbb R^n)$. At the same time, $\widetilde\varphi_j$ converges in $L^p(\mathbb R^n)$ to the extension by zero of $u$, so $\widetilde u=u$ almost everywhere on $\Omega$. The restriction of $\widetilde u$ to $\overline\Omega$ is the desired representative and, by the definition of the extension norm, satisfies \eqref{eq:morrey-primer-orden-W0}.

If two continuous representatives of the same class in $L^p(\Omega)$ agree almost everywhere, they agree throughout $\Omega$, since the set on which they differ would otherwise be open and have positive measure. Continuity up to the boundary completes uniqueness on $\overline\Omega$.
\end{proof}

The higher-order result admits a uniform formulation in terms of the Sobolev regularity number
\[
\sigma:=m-\frac{n}{p}.
\]
The restriction $0<\alpha<1$ avoids the Lipschitz endpoint, which requires separate treatment when $\sigma$ is an integer and does not in general follow from the same argument.

\begin{theorem}[Higher-order Sobolev--Morrey embedding]
\label{teo:sobolev-morrey-euclidiano}
\index{Sobolev Morrey theorem@Sobolev--Morrey theorem!higher-order}
Let $\Omega\subseteq\mathbb R^n$ be a bounded open subset, let $m\in\mathbb N$ and $1\leq p<\infty$, and suppose that
\[
\sigma:=m-\frac{n}{p}>0.
\]
If $k\in\mathbb N_0$ and $0<\alpha<1$ satisfy
\begin{equation}
\label{eq:condicion-sobolev-morrey-euclidiana}
k+\alpha\leq\sigma,
\end{equation}
then every element of $W_0^{m,p}(\Omega)$ has a unique representative in $C^{k,\alpha}(\overline\Omega)$, and there is a constant $C=C(n,m,p,k,\alpha,\Omega)>0$ such that
\begin{equation}
\label{eq:estimacion-sobolev-morrey-euclidiana}
\|u\|_{C^{k,\alpha}(\overline\Omega)}
\leq
C\|u\|_{W^{m,p}(\Omega)}.
\end{equation}

More precisely, the following assertions hold.
\begin{enumerate}[label=(\alph*)]
\item If $\sigma\notin\mathbb N$, $\ell:=\lfloor\sigma\rfloor$, and $\theta:=\sigma-\ell$, then
\[
W_0^{m,p}(\Omega)
\hookrightarrow
C^{\ell,\theta}(\overline\Omega)
\]
continuously.
\item If $\sigma=s\in\mathbb N$, then for every $0<\theta<1$,
\[
W_0^{m,p}(\Omega)
\hookrightarrow
C^{s-1,\theta}(\overline\Omega)
\]
continuously.
\end{enumerate}
\end{theorem}

\begin{proof}
For $u\in C_c^\infty(\Omega)$, let $\widetilde u$ be its zero
extension to $\mathbb R^n$. It is smooth and
\[
\|\widetilde u\|_{W^{m,p}(\mathbb R^n)}
=\|u\|_{W^{m,p}(\Omega)}.
\]
Every extension used below is a derivative of this same global function.

If $\sigma\notin\mathbb N$, set
$\ell=\lfloor\sigma\rfloor$, $\theta=\sigma-\ell$, and $d=m-\ell$.
Then $(d-1)p<n<dp$.
If $d=1$, take $q=p>n$.
If $d\geq2$, Corollary~\ref{cor:sobolev-exponente-critico-euclidiano}
gives, with $q=\displaystyle\frac{np}{n-(d-1)p}>n$,
\[
\|u\|_{W^{\ell+1,q}(\Omega)}
\leq C\|u\|_{W^{m,p}(\Omega)},
\qquad 1-\frac{n}{q}=\theta.
\]
For each $|\beta|\leq\ell$, the proof of
Theorem~\ref{teo:morrey-primer-orden-W0}, applied to
$D^\beta\widetilde u$, provides the global bound
\[
\|D^\beta\widetilde u\|_{C^{0,\theta}(\mathbb R^n)}
\leq C\|D^\beta u\|_{W^{1,q}(\Omega)}
\leq C\|u\|_{W^{m,p}(\Omega)}.
\]
The Sobolev norms of $u$ and its extension agree at this step.
Summing these bounds over the multiindices proves
\[
\|\widetilde u\|_{C^{\ell,\theta}(\mathbb R^n)}
\leq C\|u\|_{W^{m,p}(\Omega)}.
\]

If $\sigma=s\in\mathbb N$, then $d=m-s=\displaystyle\frac{n}{p}\in\mathbb N$.
For $d=1$ we already have $W_0^{m,p}(\Omega)=W_0^{s+1,n}(\Omega)$;
for $d\geq2$, the same corollary gives a continuous embedding
into $W_0^{s+1,n}(\Omega)$. Fix $0<\theta<1$ and choose $q>n$
with $\eta:=1-\displaystyle\frac{n}{q}>\theta$. The critical case of
Theorem~\ref{desigualdad de poincaré}, applied to
$D^\gamma u$ for $|\gamma|\leq s$, implies
\[
\|D^\gamma u\|_{L^q(\Omega)}
\leq C_q\|u\|_{W^{m,p}(\Omega)}.
\]
Hence, if $|\beta|\leq s-1$, the same global Morrey bound gives
\[
\|D^\beta\widetilde u\|_{C^{0,\eta}(\mathbb R^n)}
\leq C\|u\|_{W^{m,p}(\Omega)}.
\]
For distances at most one, the modulus of order $\eta$ controls
that of order $\theta$; for larger distances, use twice the uniform
norm. Therefore,
\[
\|\widetilde u\|_{C^{s-1,\theta}(\mathbb R^n)}
\leq C\|u\|_{W^{m,p}(\Omega)}.
\]

In both cases, apply the global estimate to differences of a sequence
$u_j\in C_c^\infty(\Omega)$ converging to
$u\in W_0^{m,p}(\Omega)$. Their extensions converge in the corresponding
global Hölder space, which is complete, to a single function $v$.
They also converge in $W^{m,p}(\mathbb R^n)$ to the zero extension
of $u$, by the norm identity and continuity of distributional
derivatives. Both limits agree as distributions, so the classical
derivatives of $v$ represent the weak derivatives of $u$ on $\Omega$.
The definition of the norm by extensions gives
\[
\|u\|_{C^{\ell,\theta}(\overline\Omega)}
\leq\|v\|_{C^{\ell,\theta}(\mathbb R^n)}
\]
in the noninteger case, and the same inequality with $s-1$ in the
integer case. This proves (a) and (b).

For a pair $(k,\alpha)$ satisfying
\eqref{eq:condicion-sobolev-morrey-euclidiana}, if $k=\ell$ in the
noninteger case, then $\alpha\leq\theta$ and the exponent is reduced
as above. If $k<\ell$, the order $k$ derivatives of $v$ have bounded
first derivatives and are Lipschitz; their uniform boundedness
makes them globally Hölder with any exponent in $(0,1)$.
In the integer case, $k\leq s-1$; choose $\theta>\alpha$ and use
the same argument. This gives
\eqref{eq:estimacion-sobolev-morrey-euclidiana}.
Two continuous representatives equal almost everywhere agree on
the open set and then on its closure by continuity. This also proves
uniqueness.
\end{proof}

To prove compactness, we shall need an elementary inequality between Hölder seminorms.

\begin{lemma}[Interpolation of Hölder seminorms]
\label{lem:interpolacion-seminormas-holder}
Let $K\subseteq\mathbb R^n$, $0<\alpha<\beta<1$, and $v\in C^{0,\beta}(K)$. Then
\begin{equation}
\label{eq:interpolacion-seminormas-holder}
[v]_{C^{0,\alpha}(K)}
\leq
2^{1-\frac{\alpha}{\beta}}
\|v\|_{C^0(K)}^{1-\frac{\alpha}{\beta}}
[v]_{C^{0,\beta}(K)}^{\frac{\alpha}{\beta}}.
\end{equation}
\end{lemma}

\begin{proof}
Set $M:=\|v\|_{C^0(K)}$ and $H:=[v]_{C^{0,\beta}(K)}$. For $x\neq y$, we have simultaneously
\[
|v(x)-v(y)|\leq2M,
\qquad
|v(x)-v(y)|\leq H|x-y|^\beta.
\]
If $H=0$, the assertion is immediate. If $H>0$, separate the cases $|x-y|\leq\left(\displaystyle\frac{2M}{H}\right)^{\displaystyle\frac{1}{\beta}}$ and $|x-y|>\left(\displaystyle\frac{2M}{H}\right)^{\displaystyle\frac{1}{\beta}}$. Both give
\[
\frac{|v(x)-v(y)|}{|x-y|^\alpha}
\leq
(2M)^{1-\frac{\alpha}{\beta}}H^{\frac{\alpha}{\beta}}.
\]
Taking the supremum completes the proof.
\end{proof}

\begin{lemma}[Arzelà--Ascoli on compact Euclidean subsets]
\label{lem:arzela-ascoli-compactos-euclidianos}
Let $K\subseteq\mathbb R^n$ be compact and let $(v_j)_{j\in\mathbb N}$ be a sequence of continuous functions on $K$. Suppose that there is $C>0$ such that
\[
\|v_j\|_{C^0(K)}\leq C
\qquad
\text{for every $j\in\mathbb N$},
\]
and that the family is equicontinuous. Then there is a subsequence converging uniformly on $K$.
\end{lemma}

\begin{proof}
For each $m\in\mathbb N$, compactness of $K$ provides a finite set $D_m\subseteq K$ such that the balls of radius $\displaystyle\frac{1}{m}$ centered at $D_m$ cover $K$. The set $\displaystyle D:=\bigcup_{m\in\mathbb N}D_m$ is countable and dense. Enumerate it as $D=\{x_1,x_2,\dots\}$. The numerical sequence $(v_j(x_1))_{j\in\mathbb N}$ is bounded. By Bolzano--Weierstrass, it has a convergent subsequence. Extract from this a subsequence converging at $x_2$ and continue inductively. The diagonal subsequence, again denoted by $(v_j)$, converges at each point of $D$.

We show that it is uniformly Cauchy. Let $\varepsilon>0$. By equicontinuity, there is $\delta>0$ such that
\[
|v_j(x)-v_j(y)|<\frac{\varepsilon}{3}
\]
for every $j$ and any $x,y\in K$ with $|x-y|<\delta$. Compactness of $K$ allows us to choose $x_{i_1},\dots,x_{i_N}\in D$ so that the balls of radius $\delta$ centered at these points cover $K$. Since the subsequence converges at each of these points, there is $j_0$ such that
\[
|v_j(x_{i_\nu})-v_\ell(x_{i_\nu})|
<\frac{\varepsilon}{3}
\]
for $j,\ell\geq j_0$ and $\nu\in\{1,\dots,N\}$. Given $x\in K$, choose $\nu$ with $|x-x_{i_\nu}|<\delta$. Then
\[
|v_j(x)-v_\ell(x)|
\leq
|v_j(x)-v_j(x_{i_\nu})|
+|v_j(x_{i_\nu})-v_\ell(x_{i_\nu})|
+|v_\ell(x_{i_\nu})-v_\ell(x)|
<\varepsilon.
\]
Thus, the subsequence is Cauchy in $C^0(K)$ and converges uniformly, since this space is Banach.
\end{proof}

\begin{theorem}[Rellich--Kondrashov into Hölder spaces]
\label{teo:rellich-kondrashov-holder-euclidiano}
\index{Rellich Kondrashov theorem@Rellich--Kondrashov theorem!into Hölder spaces}
Let $\Omega\subseteq\mathbb R^n$ be a bounded open subset, let $m\in\mathbb N$ and $1\leq p<\infty$, and suppose that
\[
\sigma:=m-\frac{n}{p}>0.
\]
If $k\in\mathbb N_0$ and $0<\alpha<1$ satisfy
\begin{equation}
\label{eq:condicion-compacta-holder-euclidiana}
k+\alpha<\sigma,
\end{equation}
then the embedding
\[
W_0^{m,p}(\Omega)
\hookrightarrow
C^{k,\alpha}(\overline\Omega)
\]
is compact.
\end{theorem}

\begin{proof}
Choose $\beta$ such that
\[
\alpha<\beta<1,
\qquad
k+\beta\leq\sigma.
\]
This choice is possible by \eqref{eq:condicion-compacta-holder-euclidiana}. Theorem~\ref{teo:sobolev-morrey-euclidiano} shows that every bounded sequence $(u_j)_{j\in\mathbb N}$ in $W_0^{m,p}(\Omega)$, after choosing the representatives obtained by extension by zero, is bounded in $C^{k,\beta}(\mathbb R^n)$. These extensions vanish outside $\overline\Omega$. Fix a compact cube $Q$ containing $\overline\Omega$ in its interior.

For each $|\gamma|=k$, the family $(D^\gamma u_j)$ is uniformly bounded and uniformly Hölder continuous on $Q$. If $|\gamma|<k$, the family is uniformly Lipschitz, since its first derivatives are uniformly bounded. Lemma~\ref{lem:arzela-ascoli-compactos-euclidianos} and a diagonal argument over the finite family of multi-indices allow us to extract a subsequence, which we do not relabel, such that
\[
D^\gamma u_j
\longrightarrow
v_\gamma
\qquad
\text{uniformly on $Q$}
\]
for every $|\gamma|\leq k$. The identities obtained from the fundamental theorem of calculus pass to the limit, so there is $v\in C^k(Q)$ with $D^\gamma v=v_\gamma$.

For $|\gamma|=k$, apply Lemma~\ref{lem:interpolacion-seminormas-holder} to $D^\gamma(u_j-u_\ell)$. The seminorms of order $\beta$ are uniformly bounded, while the uniform norms tend to zero. Consequently,
\[
[D^\gamma(u_j-u_\ell)]_{C^{0,\alpha}(Q)}
\longrightarrow0
\qquad
\text{as $j,\ell\to\infty$}.
\]
Thus, the subsequence is Cauchy in the $\|\cdot\|_{C^{k,\alpha}_{\mathrm{int}}(Q)}$ norm and converges to $v$ in that norm. All extensions are supported in the fixed compact set $\overline\Omega\Subset\operatorname{Int}(Q)$. Hence convergence in the $C^{k,\alpha}_{\mathrm{int}}(Q)$ norm implies convergence of these same extensions in $C^{k,\alpha}(\mathbb R^n)$: outside $Q$ the functions vanish, and pairs with one point outside $Q$ are controlled using the positive distance between $\overline\Omega$ and $\partial Q$. Restricting to $\overline\Omega$ and using the extension norm, we obtain a subsequence converging in $C^{k,\alpha}(\overline\Omega)$.
\end{proof}

\begin{remark}[Optimal form of the exponents]
\label{obs:forma-optima-sobolev-morrey-euclidiano}
Let $\sigma=m-\displaystyle\frac{n}{p}>0$.
\begin{enumerate}[label=(\alph*)]
\item If $\sigma=\ell+\theta$ with $\ell\in\mathbb N_0$ and $0<\theta<1$, the embedding
\[
W_0^{m,p}(\Omega)
\hookrightarrow
C^{\ell,\theta}(\overline\Omega)
\]
is continuous, and becomes compact upon replacing $\theta$ by any $0<\alpha<\theta$.
\item If $\sigma=s\in\mathbb N$, then for every $0<\alpha<1$ the embedding
\[
W_0^{m,p}(\Omega)
\hookrightarrow
C^{s-1,\alpha}(\overline\Omega)
\]
is compact and, in particular, continuous.
\end{enumerate}
\end{remark}

\begin{corollary}[Infinite exponents]
\label{cor:encajes-exponentes-infinitos-euclidianos}
Let $\Omega\subseteq\mathbb R^n$ be open and bounded, and let $m>k\geq0$ and $1\leq p<\infty$. If $m-\displaystyle\frac{n}{p}>k$, the embedding
\[
W_0^{m,p}(\Omega)\hookrightarrow W^{k,\infty}(\Omega)
\]
is compact. The same conclusion holds for $W^{m,p}(\Omega)$ if $\Omega$ is an extension domain for $W^{m,p}(\Omega)$. At the endpoint $p=1$ and $m-k=n$, both embeddings into $W^{k,\infty}(\Omega)$ remain continuous, although they are not compact when $\Omega$ is nonempty.

If $p=\infty$, $m>k\geq0$, and $\Omega$ is an extension domain for $W^{m,\infty}(\Omega)$, the embeddings $W^{m,\infty}(\Omega)\hookrightarrow W^{k,q}(\Omega)$ are compact for $1\leq q\leq\infty$. They are also compact on $W_0^{m,\infty}(\Omega):=
\overline{C_c^\infty(\Omega)}^{W^{m,\infty}(\Omega)}$ for any bounded open subset.

On $\mathbb R^n$, for $1\leq p<\infty$ and $m-\displaystyle\frac{n}{p}>k$, and also for $p=1$ and $m-k=n$, the embedding $W^{m,p}(\mathbb R^n)\hookrightarrow W^{k,\infty}(\mathbb R^n)$ is continuous. Restriction through an extension operator gives the same continuity on the half-space.
\end{corollary}
\begin{proof}
If $m-\displaystyle\frac{n}{p}>k$, choose $0<\alpha<1$ with $k+\alpha<m-\displaystyle\frac{n}{p}$. Theorem~\ref{teo:rellich-kondrashov-holder-euclidiano} gives compactness in $C^{k,\alpha}(\overline\Omega)$, whose norm controls all derivatives of order at most $k$ in $L^\infty(\Omega)$. On an extension domain, apply the same argument to $Tu=\chi\mathcal Eu$ on a ball, as in Theorem~\ref{rellich kondrashov generalizado}.

At the endpoint $p=1$, for $\varphi\in C_c^\infty(\mathbb R^n)$, successive applications of the fundamental theorem of calculus in the $n$ coordinates give
\[
\varphi(x)=\int_{-\infty}^{x_1}\cdots\int_{-\infty}^{x_n}
D_1\cdots D_n\varphi(t_1,\dots,t_n)\,dt_n\cdots dt_1,
\qquad
\|\varphi\|_\infty\leq\|D_1\cdots D_n\varphi\|_{L^{1}(\mathbb R^n)}.
\]
Apply this bound to each $D^\alpha\varphi$, $|\alpha|\leq k$. All derivatives on the right-hand side have order at most $k+n=m$. Density shows that the extensions by zero of the approximants converge uniformly together with their derivatives of order at most $k$. The fundamental theorem of calculus along segments identifies the limit as a function of class $C^k$ and gives continuity into $W^{k,\infty}(\Omega)$. Applying a cutoff to an extension treats $W^{m,1}(\Omega)$. To see that compactness fails, take $u_\varepsilon(x)=\varepsilon^k\varphi(\displaystyle\frac{x-x_0}{\varepsilon})$ with support in an interior ball: their $W^{m,1}(\Omega)$ norms are bounded, and some derivative of order $k$ has a constant nonzero uniform norm, although it converges to zero almost everywhere.

For $p=\infty$, fix a ball $B$ and a cutoff
$\chi\in C_c^\infty(B)$ equal to one near $\overline\Omega$.
If an extension operator exists, set
$Tu:=(\chi\mathcal Eu)\restriction_B$.
For a bounded family in $W^{m,\infty}(\Omega)$, the functions
$Tu$ form a bounded family in $W_0^{m,r}(B)$ for every
$1\leq r<\infty$, by Leibniz, Hölder, and their fixed compact
support. Choose $r$ with $m-\displaystyle\frac{n}{r}>k$ and apply the first part on
$B$. This gives a subsequence converging in $W^{k,\infty}(B)$;
restriction gives convergence in $W^{k,\infty}(\Omega)$ and,
by finite measure, in every $W^{k,q}(\Omega)$.
If $u\in W_0^{m,\infty}(\Omega)$, its test function approximants
converge in that norm. Their zero extensions converge in
$W^{m,\infty}(\mathbb R^n)$ by the norm identity and continuity
of distributional derivatives. Restricted to $B$, they also
converge in $W^{m,r}(B)$ and belong to $W_0^{m,r}(B)$.
The same estimate and the first part apply to the zero
extensions of every bounded family.

For global continuity under the strict inequality $m-\displaystyle\frac{n}{p}>k$, repeat the localization on the cubes $Q_a$ from Proposition~\ref{prop:encajes-globales-Rn-finitos}. The Sobolev--Morrey embedding, applied to $\chi_au$, controls $\|u\|_{W^{k,\infty}(Q_a)}$ by $C\|u\|_{W^{m,p}(a+(-1,2)^n)}$, with $C$ independent of $a$. Taking the supremum over $a\in\mathbb Z^n$ gives the global estimate.

In the critical case $p=1$ and $m=k+n$, the inequality obtained by successive integration gives directly, for $\varphi\in C_c^\infty(\mathbb R^n)$,
\[
\sum_{|\alpha|\leq k}\|D^\alpha\varphi\|_\infty
\leq\sum_{|\alpha|\leq k}
\|D^{\alpha+(1,\ldots,1)}\varphi\|_{L^1(\mathbb R^n)}
\leq C_{n,k}\|\varphi\|_{W^{m,1}(\mathbb R^n)}.
\]
If $\varphi_j\to u$ in $W^{m,1}(\mathbb R^n)$, applying this bound to $\varphi_j-\varphi_i$ makes all derivatives of order at most $k$ converge uniformly. Their limits are the distributional derivatives of $u$: for each test function, we pass to the limit using the $L^1$ convergence of $\varphi_j$ and the uniform convergence of $D^\alpha\varphi_j$. Density therefore gives the same inequality for $u$ and proves the global critical embedding. In both cases, composing with the extension operator and restricting to the half-space gives the asserted continuity there.
\end{proof}

Equality of indices with an infinite exponent requires its own analysis. For example, in dimension one, $W_0^{1,1}(a,b)\hookrightarrow C^0([a,b])$ continuously: the estimate $\|\varphi\|_\infty\leq\|\varphi'\|_{L^{1}((a,b))}$ and density give the continuous representative. This embedding is not compact. If $c\in(a,b)$ and $0\leq\varphi\in C_c^\infty((-1,1))$ satisfies $\varphi(0)=1$, the functions $u_j(x)=\varphi(j(x-c))$, for sufficiently large $j$, are bounded in $W_0^{1,1}(a,b)$, have uniform norm one, and converge to zero at every point other than $c$. They have no uniformly convergent subsequence. The preceding assertions do not infer compactness at any equality endpoint solely from a comparison of indices.

\section{Localized embeddings on the half-space}

The upper half-space is the Euclidean model for boundary charts. Since it is unbounded, compactness cannot be formulated globally without a condition preventing functions from moving toward infinity. In applications to compact manifolds, each function localized by a partition of unity has support contained in a fixed compact set. This is exactly the situation recorded below.

\begin{definition}[Uniformly Lipschitz continuous boundary]
\label{def:frontera-uniformemente-lipschitz}
\index{boundary!uniformly Lipschitz}
Let $\Omega\subseteq\mathbb R^n$ be an open subset. We say that its boundary $\partial\Omega$ is \textbf{uniformly Lipschitz continuous} if there exist $\varepsilon,L>0$, $M\in\mathbb N$, and a countable locally finite open cover $(\Omega_j)_{j\in\mathbb N}$ of $\partial\Omega$ such that:
\begin{enumerate}
\item for each $x\in\partial\Omega$, there is $j\in\mathbb N$ such that $B_{\mathrm{euc}}(x,\varepsilon)\subseteq\Omega_j$;
\item no point of $\mathbb R^n$ belongs to more than $M$ sets in the family $(\Omega_j)_{j\in\mathbb N}$;
\item for each $j\in\mathbb N$, there exist a coordinate system $y=(y',y_n)\in\mathbb R^{n-1}\times\mathbb R$ and a Lipschitz function $f_j\colon \mathbb R^{n-1}\longrightarrow\mathbb R$, both depending on $j$, such that $\operatorname{Lip}(f_j)\leq L$ and
\[
\Omega_j\cap\Omega=\Omega_j\cap V_j,
\qquad
V_j:=\{(y',y_n)\in\mathbb R^{n-1}\times\mathbb R\mid y_n>f_j(y')\}.
\]
\end{enumerate}
The numbers $\varepsilon$, $L$, and $M$ are called the \textbf{uniformity parameters} of the boundary.
\end{definition}

\begin{theorem}[Extension on open subsets with uniformly Lipschitz boundary]
\label{teo:extension-sobolev-uniforme-lipschitz}
\index{extension operator!on open subsets with uniformly Lipschitz boundary}
Let $\Omega\subseteq\mathbb R^n$ be an open subset whose boundary is uniformly Lipschitz continuous with parameters $\varepsilon,L>0$ and $M\in\mathbb N$. Let $1\leq p\leq\infty$ and $m\in\mathbb N$. Then there is a continuous linear operator
\[
\mathcal E\colon W^{m,p}(\Omega)\longrightarrow W^{m,p}(\mathbb R^n)
\]
such that, for every $u\in W^{m,p}(\Omega)$, we have $\mathcal E(u)(x)=u(x)$ for almost every $x\in\Omega$. Moreover, there is a constant $c=c(n,m,p)>0$ for which
\begin{equation}
\label{eq:estimacion-extension-uniforme-Lp}
\|\mathcal E(u)\|_{L^p(\mathbb R^n)}
\leq
c(1+M)\|u\|_{L^p(\Omega)}
\end{equation}
and, for each $k\in\{1,\dots,m\}$,
\begin{equation}
\label{eq:estimacion-extension-uniforme-derivadas}
\|\nabla^k\mathcal E(u)\|_{L^p(T^{(0,k)}(\mathbb R^n))}
\leq
c\bigl(1+M(1+L^k)\bigr)
\displaystyle\sum_{i=0}^{k}
\left(\frac{M}{\varepsilon}\right)^{k-i}
\|\nabla^iu\|_{L^p(T^{(0,i)}(\Omega))},
\end{equation}
where $\nabla^0u:=u$ and the covariant derivatives are those induced by the Euclidean metric. In particular, there is a constant $C=C(n,m,p,\varepsilon,M,L)>0$ such that
\[
\|\mathcal E(u)\|_{W^{m,p}(\mathbb R^n)}
\leq
C\|u\|_{W^{m,p}(\Omega)}.
\]
\end{theorem}

This is Stein's extension theorem. The construction of the operator and the estimates \eqref{eq:estimacion-extension-uniforme-Lp} and \eqref{eq:estimacion-extension-uniforme-derivadas} can be found in \cite[Theorem~13.17]{Leoni2017}.

\begin{corollary}[Embeddings on a ball without boundary conditions]
\label{cor:encajes-bola-espacio-completo}
Let $B\subset\mathbb R^n$ be a bounded open ball. Let $\ell>k\geq0$ be integers and $1\leq p,q<\infty$. If
\[
\ell-\frac{n}{p}\geq k-\frac{n}{q},
\]
then $W^{\ell,p}(B)\hookrightarrow W^{k,q}(B)$ continuously; if the inequality is strict, the embedding is compact.
\end{corollary}
\begin{proof}
The boundary of a ball is smooth. Near each point, after a translation and a rotation, it can be represented as the graph of a smooth function on a ball in $\mathbb R^{n-1}$. Shrinking this ball makes the derivative of the function bounded, and the mean value theorem then shows that the graph is Lipschitz. Multiplying the function by a cutoff equal to one on a smaller ball gives a Lipschitz function defined on all of $\mathbb R^{n-1}$ without changing this local representation. Compactness of $\partial B$ allows us to choose finitely many representations on open sets $V_i$ and smaller open sets $U_i$ that still cover $\partial B$, with $\overline{U_i}\Subset V_i$. The minimum of the positive distances $\operatorname{dist}(\overline{U_i},\mathbb R^n\setminus V_i)$, the maximum of the Lipschitz constants, and the number of charts give the parameters in Definition~\ref{def:frontera-uniformemente-lipschitz}. For $n=1$, the boundary consists of the two endpoints of an interval, and the same description reduces to half-lines.

By Theorem~\ref{teo:extension-sobolev-uniforme-lipschitz}, $B$ is an extension domain for $W^{\ell,p}(B)$. Apply the second part of Theorem~\ref{prop:encaje-continuo-indices-finitos} and, under the strict inequality, of Theorem~\ref{rellich kondrashov generalizado}. Both results apply to the full space $W^{\ell,p}(B)$, without requiring interior support or membership in $W_0^{\ell,p}(B)$.
\end{proof}

\begin{corollary}[Half-space extension operator]
\label{teo:extension-sobolev-semiespacio}
\index{extension operator!half-space}
Let $k\in\mathbb N_0$ and $1\leq p\leq\infty$. There is a continuous linear operator
\[
\mathcal E_k\colon W^{k,p}(\mathbb R^n_+)\longrightarrow W^{k,p}(\mathbb R^n)
\]
such that $(\mathcal E_ku)\restriction_{\mathbb R^n_+}=u$ at almost every point of $\mathbb R^n_+$. In particular, there is $C_{k,p,n}>0$ such that
\[
\|\mathcal E_ku\|_{W^{k,p}(\mathbb R^n)}
\leq
C_{k,p,n}\|u\|_{W^{k,p}(\mathbb R^n_+)}
\]
for every $u\in W^{k,p}(\mathbb R^n_+)$.
\end{corollary}

\begin{proof}
If $k=0$, extension by zero suffices. If $k\geq1$, the half-space
$\mathbb R^n_+:=\{(x',t)\in\mathbb R^{n-1}\times\mathbb R\mid t>0\}$
has uniformly Lipschitz boundary, and Theorem~\ref{teo:extension-sobolev-uniforme-lipschitz} gives the operator for every $1\leq p\leq\infty$.

For finite exponents, one can see explicitly how higher-order derivatives are preserved. Choose $a_1,\dots,a_k\in\mathbb R$ solving
\[
\sum_{j=1}^{k}a_j(-j)^r=1,
\qquad r\in\{0,\dots,k-1\}.
\]
The matrix is a Vandermonde matrix at the distinct points $-1,\dots,-k$, and is therefore invertible. For a function smooth up to the boundary, define
\[
(\mathcal E_ku)(x',t):=
\begin{cases}
u(x',t),&t>0,\\
\displaystyle\sum_{j=1}^{k}a_ju(x',-jt),&t<0.
\end{cases}
\]
The tangential and normal derivatives of total order at most $k-1$ agree on both sides of $t=0$ by the preceding equations. Integration by parts on both half-spaces cancels the boundary terms. Applying this argument successively to each derivative of order at most $k-1$ shows that the weak derivatives of order at most $k$ are the ordinary derivatives on each side. For $\alpha=(\alpha',r)$, $|\alpha|\leq k$, the change of variable $s=-jt$ and Minkowski's inequality give
\[
\|D^\alpha\mathcal E_ku\|_{L^p(\{t<0\})}
\leq\sum_{j=1}^{k}|a_j|j^{r-\frac{1}{p}}
\|D^\alpha u\|_{L^p(\mathbb R^n_+)}.
\]
The extension operator already obtained and Proposition~\ref{prop:densidad-sobolev-orden-entero} show that restrictions of $C_c^\infty(\mathbb R^n)$ are dense in $W^{k,p}(\mathbb R^n_+)$ for $p<\infty$. The estimate allows this formula to be extended by continuity and verifies all its derivatives. For $k=1$, we obtain even reflection; for $k\geq2$, the finite combination also ensures agreement of the necessary normal derivatives.
\end{proof}

\begin{remark}[Density and boundary conditions on the half-space]
\label{obs:densidad-trazas-semiespacio-euclidiano}
Restrictions of $C_c^\infty(\mathbb R^n)$ are dense in $W^{m,p}(\mathbb R^n_+)$ for $m\geq1$ and $1\leq p<\infty$. These functions need not vanish on the boundary. In contrast, $W_0^{m,p}(\mathbb R^n_+)$ is defined using approximants compactly supported in the open half-space. Theorem~\ref{teo:traza-semiespacio} will characterize this latter space by the vanishing of traces of order less than $m$.
\end{remark}

For a compact set $K\subseteq\overline{\mathbb R^n_+}$, we shall write
\[
W_K^{k,p}(\mathbb R^n_+)
:=
\left\{u\in W^{k,p}(\mathbb R^n_+)\middle|\operatorname{supp}u\subseteq K\right\},
\]
where support is understood in the essential sense within $\overline{\mathbb R^n_+}$.

\begin{theorem}[Localized Sobolev and Rellich--Kondrashov embeddings on the half-space]
\label{teo:encajes-localizados-semiespacio}
\index{Sobolev embedding theorem!on the half-space}
\index{Rellich Kondrashov theorem@Rellich--Kondrashov theorem!on the half-space}
Let $K\subseteq\overline{\mathbb R^n_+}$ be compact, let $m>k\geq0$ be integers, and let $1\leq p,q<\infty$. If
\[
m-\frac{n}{p}\geq k-\frac{n}{q},
\]
the inclusion
$W_K^{m,p}(\mathbb R^n_+)\hookrightarrow W^{k,q}(\mathbb R^n_+)$
is continuous. If the inequality is strict, this inclusion is compact. In particular, $W_K^{m,p}(\mathbb R^n_+)\hookrightarrow
W^{m-1,p}(\mathbb R^n_+)$ is compact for every $m\geq1$ and $1\leq p<\infty$.

If $q=\infty$ and $p<\infty$, the inclusion is compact when $m-\displaystyle\frac{n}{p}>k$; for $p=1$ and $m-k=n$, it is continuous. If $p=\infty$, it is compact for every $1\leq q\leq\infty$ and $m>k$. These assertions also hold for the subspace $W_K^{m,p}(\mathbb R^n_+)\cap W_0^{m,p}(\mathbb R^n_+)$; for finite exponents, its limit belongs to $W_0^{k,q}(\mathbb R^n_+)$.
\end{theorem}
\begin{proof}
Take $R>0$ such that $K\subseteq B_{\mathrm{euc}}(0,R)\cap\overline{\mathbb R^n_+}$, and $\chi\in C_c^\infty(B_{\mathrm{euc}}(0,2R))$ equal to one on a neighborhood of $K$. Define
\[
Tu:=\chi\mathcal E_mu.
\]
The Leibniz rule and Corollary~\ref{teo:extension-sobolev-semiespacio} give
\[
\|Tu\|_{W^{m,p}(\mathbb R^n)}
\leq C\|u\|_{W^{m,p}(\mathbb R^n_+)}.
\]
Moreover, $Tu$ has compact support contained in $B:=B_{\mathrm{euc}}(0,2R)$, and $(Tu)\restriction_{\mathbb R^n_+}=u$, since $u$ vanishes outside $K$. For $p<\infty$, Proposition~\ref{prop:extension-cero-W0-orden-entero} gives $Tu\in W_0^{m,p}(B)$. Continuity for finite exponents follows from Theorem~\ref{prop:encaje-continuo-indices-finitos} on $B$, and compactness from Theorem~\ref{rellich kondrashov generalizado}. In the compactness assertion, the resulting sequence converges to an element of $W_0^{k,q}(B)$, since this is the target space in Theorem~\ref{rellich kondrashov generalizado}. Extension by zero is continuous by Proposition~\ref{prop:extension-cero-W0-orden-entero}. Hence this convergence is preserved in $W^{k,q}(\mathbb R^n)$ and, upon restriction, in $W^{k,q}(\mathbb R^n_+)$. The same cutoff and the same ball are used for the entire sequence.

For $p<\infty$ and $q=\infty$, use Corollary~\ref{cor:encajes-exponentes-infinitos-euclidianos} on $B$. For $p=\infty$, the functions $Tu$ are bounded in $W^{m,r}(\mathbb R^n)$ for each $r<\infty$ because their support is fixed. Choose $r>\displaystyle\frac{n}{m-k}$ and apply that corollary again; the finite exponents follow from the finite measure of the support.

If each $u_j$ also belongs to $W_0^{m,p}(\mathbb R^n_+)$ and the exponents are finite, we may use $Z_{\mathbb R^n_+}u_j$ in place of $\mathcal E_mu_j$. To justify the zero condition in the target space, approximate each $u_j$ in $W^{m,p}(\mathbb R^n_+)$ by test functions and multiply them by $\chi$. Multiplication preserves the approximation and places all supports in $B$. The continuous embedding on $B$ shows that the restrictions converge in $W^{k,q}(\mathbb R^n_+)$ to $u_j$, so $u_j\in W_0^{k,q}(\mathbb R^n_+)$. Closedness of this subspace preserves the condition in the limit supplied by compactness.
\end{proof}

\begin{remark}[Support up to the boundary and local compactness]
\label{obs:soportes-semiespacio-frontera}
The compact set $K$ in the theorem may meet $\{t=0\}$. In this case, a function in $W_K^{m,p}(\mathbb R^n_+)$ may have nonzero trace, and its extension by zero may fail to belong to $W^{m,p}(\mathbb R^n)$; this is why $\mathcal E_m$ is used. If $K\Subset\mathbb R^n_+$ is contained in the interior, choose a cutoff in $C_c^\infty(\mathbb R^n_+)$ equal to one near $K$. The weak Leibniz rule and interior mollification then show that $W_K^{m,p}(\mathbb R^n_+)\subseteq W_0^{m,p}(\mathbb R^n_+)$ for $p<\infty$.

For a bounded sequence in $W^{m,p}(\mathbb R^n_+)$, cutoffs on balls and a diagonal extraction give local compactness up to the boundary: convergence holds in $W^{k,q}(B_{\mathrm{euc}}(0,R)\cap\mathbb R^n_+)$ for every $R>0$ in the strict finite range. Without a condition on the tails, the whole half-space lacks global compactness. In dimension $n\geq2$, tangential translations of a nonzero test function have constant norms and may have disjoint supports. In dimension $n=1$, translate that function toward $+\infty$. Both examples reproduce the separation of norms in Remark~\ref{obs:no-compacidad-global-Rn}, even within $W_0^{m,p}(\mathbb R^n_+)$.
\end{remark}

To apply the half-space embeddings to functions defined on a half-ball, we must first extend them across its spherical portion. Separation of the support from that portion prevents jumps; no vanishing on the flat boundary is required.

\begin{lemma}[Extension by zero from a half-ball]\label{lem:extension-cero-media-bola-semiespacio}\index{extension by zero from a half-ball}
Let $m\in\mathbb N_0$, $1\leq p<\infty$, and $R>0$. Let \(B^+=B_{\mathrm{euc}}(0,R)\cap\mathbb R^n_+\) and \(B'=B_{\mathrm{euc}}(0,R)\cap\partial\mathbb R^n_+\). Suppose that \(v\in W^{m,p}(B^+)\) and that there is a compact set \(K\subset B_{\mathrm{euc}}(0,R)\cap\{x_n\geq0\}\), lying at a positive distance from the spherical portion \(\partial B_{\mathrm{euc}}(0,R)\cap\{x_n\geq0\}\), such that \(v=0\) almost everywhere on \(B^+\setminus K\).

Define \(\widetilde v\colon \mathbb R^n_+\longrightarrow\mathbb R\) by
\[
\widetilde v(x)
=
\begin{cases}
v(x), & x\in B^+,\\
0, & x\in\mathbb R^n_+\setminus B^+.
\end{cases}
\]
Then \(\widetilde v\in W^{m,p}(\mathbb R^n_+)\), and for every multi-index \(\beta\) with \(|\beta|\leq m\), we have \(D^\beta\widetilde v=\widetilde{D^\beta v}\) in the weak sense, where \(\widetilde{D^\beta v}\) denotes the extension by zero of \(D^\beta v\) to \(\mathbb R^n_+\). In particular,
\[
\|\widetilde v\|_{W^{m,p}(\mathbb R^n_+)}
=
\|v\|_{W^{m,p}(B^+)}.
\]
\end{lemma}
\begin{figura}[htbp]
    \includegraphics[width=0.9\linewidth]{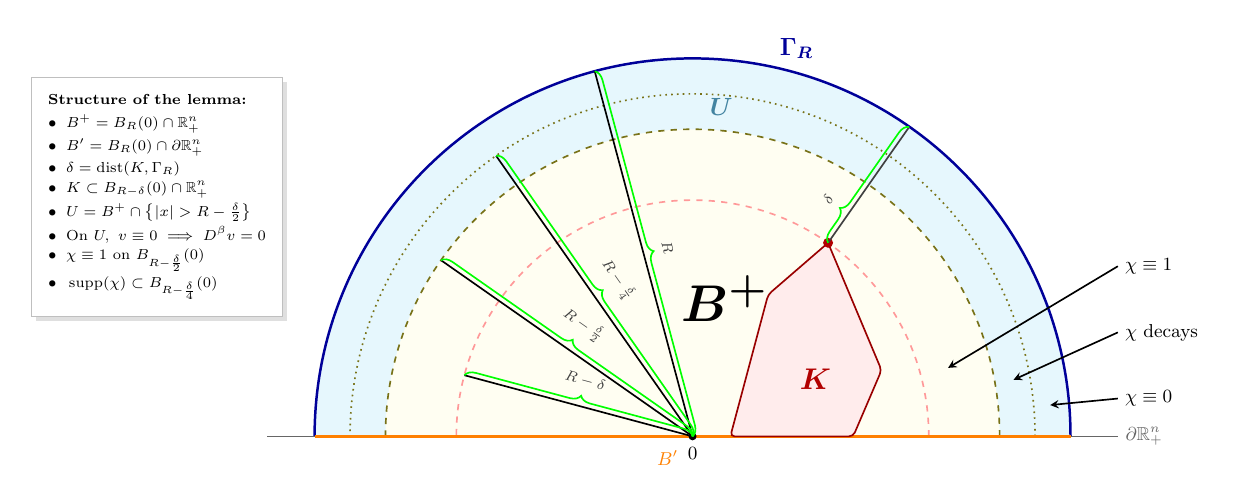}
    \caption{Outline of the proof of extension by zero. The essential support of \(v\) is contained in a compact set \(K\) separated from the spherical portion \(\Gamma_R\) by a positive distance \(\delta\). This separation allows us to construct a neighborhood \(U\) of \(\Gamma_R\) on which \(v\) and its weak derivatives vanish, as well as a cutoff function \(\chi\) equal to \(1\) near \(K\) and vanishing before it reaches \(\Gamma_R\). This prevents boundary terms from appearing on the spherical portion and gives \(D^\beta\widetilde{v}=\widetilde{D^\beta v}\).}
    \label{fig:lema-extension-por-cero-semiespacio}
\end{figura}

\begin{proof}
If \(m=0\), the definition of extension by zero gives \(\|\widetilde v\|_{L^p(\mathbb R^n_+)}=\|v\|_{L^p(B^+)}\), and the identity corresponding to \(\beta=0\) is immediate. Since \(W^{0,p}(B^+)=L^p(B^+)\), this proves the result in that case. Henceforth assume \(m\geq1\).

Denote the spherical portion of the half-ball by \(\Gamma_R:=\partial B_{\mathrm{euc}}(0,R)\cap\{x_n\geq0\}\). If \(K=\varnothing\), then \(v=0\) almost everywhere on \(B^+\), and the result is immediate. Thus assume \(K\neq\varnothing\).

Since \(K\) lies at a positive distance from \(\Gamma_R\), we may take \(\delta:=\operatorname{dist}(K,\Gamma_R)>0\). First consider what this implies. For each \(x\in K\), we have \(\operatorname{dist}(x,\Gamma_R)=R-\|x\|\). If \(x\neq0\), the point \(\displaystyle\frac{R}{\|x\|}x\) belongs to \(\Gamma_R\), and hence \(\operatorname{dist}(x,\Gamma_R)\leq\left\|x-\displaystyle\frac{R}{\|x\|}x\right\|=R-\|x\|\). The reverse inequality follows from \(\|x-\mathbf{y}\|\geq\bigl|\|\mathbf{y}\|-\|x\|\bigr|=R-\|x\|\) for every \(\mathbf{y}\in\Gamma_R\). If \(x=0\), the equality is also clear. Therefore, \(R-\|x\|\geq\delta\) for every \(x\in K\); that is, \(\|x\|\leq R-\delta\) for every \(x\in K\). In particular, \(K\subset B_{\mathrm{euc}}(0,R-\displaystyle\frac{\delta}{2})\cap\{x_n\geq0\}\).

Define \(U:=B^+\cap\left\{x\in\mathbb R^n\mid \|x\|>R-\displaystyle\frac{\delta}{2}\right\}\). Then \(U\) is an interior strip adjacent to \(\Gamma_R\), and moreover \(U\cap K=\varnothing\). Since \(v=0\) almost everywhere on \(B^+\setminus K\), it follows that \(v=0\) almost everywhere on \(U\).

By locality of weak derivatives, \(D^\beta v=0\) almost everywhere on \(U\) for every \(|\beta|\leq m\). If \(w\in W^{k,p}(B^+)\) vanishes almost everywhere on a relatively open subset \(V\subset B^+\), then for every \(\psi\in C_c^\infty(V)\), viewed as a test function on \(B^+\), we have
\[
\int_V D^\beta w\,\psi\,dx
=
\int_{B^+}D^\beta w\,\psi\,dx
=
(-1)^{|\beta|}\int_{B^+}w\,D^\beta\psi\,dx
=
(-1)^{|\beta|}\int_V w\,D^\beta\psi\,dx
=
0.
\]
For each $|\beta|\leq k$, since \(D^\beta w\in L^p(V)\), it follows that \(D^\beta w=0\) almost everywhere on \(V\).

We now prove that \(D^\beta\widetilde v=\widetilde{D^\beta v}\) in the weak sense for every \(|\beta|\leq m\). Fix such an \(\beta\) and take \(\varphi\in C_c^\infty(\mathbb R^n_+)\). Choose a function \(\chi\in C_c^\infty(\mathbb R^n)\) such that \(\chi\equiv1\) on \(B_{\mathrm{euc}}(0,R-\displaystyle\frac{\delta}{2})\) and \(\operatorname{supp}\chi\subset B_{\mathrm{euc}}(0,R-\displaystyle\frac{\delta}{4})\). This choice is possible because \(R-\displaystyle\frac{\delta}{2}<R-\displaystyle\frac{\delta}{4}<R\). In particular, \(\chi\equiv1\) on a neighborhood of \(K\), and \((\chi|_{\mathbb R^n_+})\varphi\restriction_{B^+}\in C_c^\infty(B^+)\), since \(\chi\) vanishes near the spherical portion and \(\varphi\) has compact support in \(\mathbb R^n_+\).

Since \(\widetilde v=0\) on \(\mathbb R^n_+\setminus B^+\), we have
\[
\int_{\mathbb R^n_+}\widetilde v\,D^\beta\varphi\,dx
=
\int_{B^+}v\,D^\beta\varphi\,dx.
\]
Moreover, since \(v=0\) almost everywhere on \(B^+\setminus K\), it suffices to compare the integrands on a neighborhood of \(K\). On that neighborhood, \(\chi\equiv1\), and therefore \(D^\beta(\chi\varphi)=D^\beta\varphi\), where \(\chi\varphi\) denotes \((\chi|_{\mathbb R^n_+})\varphi\). It follows that \(vD^\beta\varphi=vD^\beta(\chi\varphi)\) almost everywhere on \(B^+\). Thus,
\[
\int_{B^+}v\,D^\beta\varphi\,dx
=
\int_{B^+}v\,D^\beta(\chi\varphi)\,dx.
\]
Since \((\chi|_{\mathbb R^n_+})\varphi\restriction_{B^+}\in C_c^\infty(B^+)\), the definition of weak derivative gives
\[
\int_{B^+}v\,D^\beta(\chi\varphi)\,dx
=
(-1)^{|\beta|}\int_{B^+}D^\beta v\,\chi\varphi\,dx.
\]
Now observe that within \(B^+\), the set where \(\chi\neq1\) is contained in \(U\), except possibly on the sphere \(\left\{\|x\|=R-\displaystyle\frac{\delta}{2}\right\}\), which has measure zero. Since we already showed that \(D^\beta v=0\) almost everywhere on \(U\), we obtain
\[
\int_{B^+}D^\beta v\,\chi\varphi\,dx
=
\int_{B^+}D^\beta v\,\varphi\,dx
=
\int_{\mathbb R^n_+}\widetilde{D^\beta v}\,\varphi\,dx.
\]
Combining the preceding equalities gives
\[
\int_{\mathbb R^n_+}\widetilde v\,D^\beta\varphi\,dx
=
(-1)^{|\beta|}
\int_{\mathbb R^n_+}\widetilde{D^\beta v}\,\varphi\,dx.
\]
This proves that \(D^\beta\widetilde v=\widetilde{D^\beta v}\) in the weak sense for every \(|\beta|\leq m\).

Since \(\widetilde{D^\beta v}\in L^p(\mathbb R^n_+)\) for every \(|\beta|\leq m\), we conclude that \(\widetilde v\in W^{m,p}(\mathbb R^n_+)\). Moreover,
\[
\|\widetilde v\|_{W^{m,p}(\mathbb R^n_+)}^p
=
\displaystyle\sum_{|\beta|\leq m}\int_{\mathbb R^n_+}|D^\beta\widetilde v|^p\,dx
=
\displaystyle\sum_{|\beta|\leq m}\int_{\mathbb R^n_+}|\widetilde{D^\beta v}|^p\,dx
=
\displaystyle\sum_{|\beta|\leq m}\int_{B^+}|D^\beta v|^p\,dx
=
\|v\|_{W^{m,p}(B^+)}^p.
\]
Therefore, \(\|\widetilde v\|_{W^{m,p}(\mathbb R^n_+)}=\|v\|_{W^{m,p}(B^+)}\).
\end{proof}

\begin{corollary}[Approximation and localized embeddings on a half-ball]
\label{cor:aproximacion-encajes-media-bola}
\index{Sobolev embeddings!localized on a half-ball}
Let $B^+=B_{\mathrm{euc}}(0,R)\cap\mathbb R^n_+$, $B'=B_{\mathrm{euc}}(0,R)\cap\{x_n=0\}$, and let $K\subset B_{\mathrm{euc}}(0,R)\cap\{x_n\geq0\}$ be compact.
\begin{enumerate}
\item If $s\in\mathbb N_0$, $1\leq p<\infty$, and $v\in W^{s,p}(B^+)$ vanishes almost everywhere on $B^+\setminus K$, there exist $f_\nu\in C_c^\infty(B_{\mathrm{euc}}(0,R))$ whose restrictions converge to $v$ in $W^{s,p}(B^+)$. Their supports are contained in a common compact subset of $B_{\mathrm{euc}}(0,R)$.
\item If $\ell>k\geq0$ are integers, $1\leq p,q<\infty$, and $\ell-\displaystyle\frac{n}{p}\geq k-\displaystyle\frac{n}{q}$, there is $C>0$ such that every $v\in W^{\ell,p}(B^+)$ with that support satisfies
\[
\|v\|_{W^{k,q}(B^+)}\leq C\|v\|_{W^{\ell,p}(B^+)}.
\]
If the index inequality is strict, every bounded sequence of these functions, with the same $K$, has a subsequence converging in $W^{k,q}(B^+)$.
\end{enumerate}
The restrictions in the first part are smooth up to $B'$, but need not vanish there.
\end{corollary}
\begin{proof}
Extend $v$ by zero within the half-space using Lemma~\ref{lem:extension-cero-media-bola-semiespacio}. The function $\widetilde v$ has the same Sobolev norm as $v$ and vanishes outside $K$. For the first part, restrictions of $C_c^\infty(\mathbb R^n)$ are dense in $W^{s,p}(\mathbb R^n_+)$: for $s\geq1$, apply the operator from Corollary~\ref{teo:extension-sobolev-semiespacio} and Proposition~\ref{prop:densidad-sobolev-orden-entero}; for $s=0$, approximate the extension by zero in $L^p(\mathbb R^n)$. Thus, there exist $h_\nu\in C_c^\infty(\mathbb R^n)$ with $h_\nu\restriction_{\mathbb R^n_+}\to\widetilde v$ in $W^{s,p}(\mathbb R^n_+)$.

Choose a fixed cutoff $\eta\in C_c^\infty(B_{\mathrm{euc}}(0,R))$ equal to one near $K$ and set $f_\nu:=\eta h_\nu$. Since $\eta\widetilde v=\widetilde v$, the Leibniz rule and restriction give
\[
\begin{aligned}
\|f_\nu\restriction_{B^+}-v\|_{W^{s,p}(B^+)}
&\leq\|\eta(h_\nu\restriction_{\mathbb R^n_+}-\widetilde v)
      \|_{W^{s,p}(\mathbb R^n_+)}\\
&\leq C_{s,p,\eta}
\|h_\nu\restriction_{\mathbb R^n_+}-\widetilde v
      \|_{W^{s,p}(\mathbb R^n_+)}\longrightarrow0.
\end{aligned}
\]
The support of each $f_\nu$ is contained in $\operatorname{supp}\eta\Subset B_{\mathrm{euc}}(0,R)$.

For the second part, let $B:=B_{\mathrm{euc}}(0,2R)$ and define
\[
z:=(\mathcal E_\ell\widetilde v)\restriction_B.
\]
Half-space extension and Lemma~\ref{lem:restriccion-extension-local-euclidiana} imply $z\restriction_{B^+}=v$ and
\[
\|z\|_{W^{\ell,p}(B)}
\leq\|\mathcal E_\ell\widetilde v\|_{W^{\ell,p}(\mathbb R^n)}
\leq C_{\mathcal E}\|v\|_{W^{\ell,p}(B^+)}.
\]
By Corollary~\ref{cor:encajes-bola-espacio-completo},
\[
\|v\|_{W^{k,q}(B^+)}
\leq\|z\|_{W^{k,q}(B)}
\leq C_BC_{\mathcal E}\|v\|_{W^{\ell,p}(B^+)}.
\]
In this application we do not need to apply a cutoff to $\mathcal E_\ell\widetilde v$: the embedding is applied to its restriction in the full space on the ball. The cutoff in the first part serves a different purpose, keeping the supports of the approximants within the chart.

If $(v_\nu)_{\nu\in\mathbb N}$ is bounded and the inequality is strict, use the same operator and the same ball to define $z_\nu$. Then
\[
\sup_{\nu\in\mathbb N}\|z_\nu\|_{W^{\ell,p}(B)}
\leq C_{\mathcal E}\sup_{\nu\in\mathbb N}
\|v_\nu\|_{W^{\ell,p}(B^+)}<\infty.
\]
The same corollary provides a subsequence $z_{\nu_h}\to z_\infty$ in $W^{k,q}(B)$. Finally,
\[
\|v_{\nu_h}-z_\infty\restriction_{B^+}\|_{W^{k,q}(B^+)}
\leq\|z_{\nu_h}-z_\infty\|_{W^{k,q}(B)}\longrightarrow0.
\]
This proves compactness without imposing a zero boundary condition.
\end{proof}

\subsection{Summary of domains and embeddings}

In the table, $m>k\geq0$ are integers, $p,q$ are finite unless otherwise stated, and $\delta:=(m-\displaystyle\frac{n}{p})-(k-\displaystyle\frac{n}{q})$. The extension condition concerns the source space. On every open subset, $W_0^{m,p}(\Omega)$ is defined by closure; interior smooth density and the equality $W_0=W$ on $\mathbb R^n$ are given in Proposition~\ref{prop:densidad-sobolev-orden-entero}.

\begin{center}
\small
\begin{tabular}{p{0.19\textwidth}p{0.22\textwidth}p{0.22\textwidth}p{0.25\textwidth}}
\hline
Domain and space & Assumptions and range & Conclusion & Reference\\
\hline
$\Omega$ bounded, $W_0^{m,p}(\Omega)$ & $1\leq p,q<\infty$ &
Continuous if $\delta\geq0$; compact if $\delta>0$, into $W_0^{k,q}(\Omega)$ &
\ref{prop:encaje-continuo-indices-finitos}, \ref{rellich kondrashov generalizado}\\[2mm]
$\Omega$ bounded, $W^{m,p}(\Omega)$ & Extension domain; $1\leq p,q<\infty$ &
Continuous if $\delta\geq0$; compact if $\delta>0$ &
\ref{prop:encaje-continuo-indices-finitos}, \ref{rellich kondrashov generalizado}\\[2mm]
$\mathbb R^n$, $W^{m,p}(\mathbb R^n)=W_0^{m,p}(\mathbb R^n)$ & $1\leq p\leq q<\infty$ &
Globally continuous if $\delta\geq0$ &
\ref{prop:encajes-globales-Rn-finitos}\\[2mm]
$\mathbb R^n$, support in compact $K$ & $1\leq p,q<\infty$; $\delta>0$ &
Compact with common support; locally compact &
\ref{cor:rellich-local-Rn}\\[2mm]
$\mathbb R^n_+$, $W_K^{m,p}(\mathbb R^n_+)$; also $W_K^{m,p}(\mathbb R^n_+)\cap W_0^{m,p}(\mathbb R^n_+)$ &
$K\subseteq\overline{\mathbb R^n_+}$ compact; $1\leq p,q<\infty$ &
Continuous if $\delta\geq0$; compact if $\delta>0$ &
\ref{teo:encajes-localizados-semiespacio}\\[2mm]
$\Omega$ bounded; $W_0^{m,p}(\Omega)$ or $W^{m,p}(\Omega)$ with extension &
$p<\infty$, $q=\infty$, $m-\frac{n}{p}>k$; or $p=\infty$ &
Compact into $W^{k,q}(\Omega)$ in the indicated ranges &
\ref{cor:encajes-exponentes-infinitos-euclidianos}\\
\hline
\end{tabular}
\end{center}

\section{The trace theorem on the half-space}
\label{sec:traza-entera-euclidiana}

The trace assigns boundary values to a Sobolev function through smooth approximation. On the half-space, an estimate in the normal direction allows its construction for every $1\leq p<\infty$; applying it to derivatives will also give the characterization of $W_0^{m,p}(\mathbb R^n_+)$.

We shall use the convention in Definition~\ref{def:espacios-holder-euclidianos}: a function belongs to $C^\infty(\overline\Omega)$ if it is the restriction to $\overline\Omega$ of a smooth function defined on an open neighborhood of $\overline\Omega$.
Moreover,
\[
C_c^\infty(\overline\Omega)
=
\left\{
 f\in C^\infty(\overline\Omega)
 \middle|
 \operatorname{supp}_{\overline\Omega}(f)\text{ is compact}
\right\}.
\]

In particular, on the half-space we shall work with $C_c^\infty(\overline{\mathbb R^n_+})$. For each $m\in\mathbb N$ and $1\leq p<\infty$, this space is dense in $W^{m,p}(\mathbb R^n_+)$, and ordinary restriction to the flat boundary initially defines an operator on these functions. The trace theorem states that this operator extends uniquely and continuously to all of $W^{m,p}(\mathbb R^n_+)$.

\begin{theorem}[Trace theorem on the half-space]\label{teo:traza-semiespacio}\index{trace theorem on the half-space}
Let \(\mathbb R^n_+=\{x\in\mathbb R^n\mid x_n>0\}\), let \(\partial\mathbb R^n_+=\{x_n=0\}\cong\mathbb R^{n-1}\), and let \(m\geq1\) and \(1\leq p<\infty\). Then, for each multi-index \(\beta\) with \(|\beta|\leq m-1\), the operator
\[
f\longmapsto (D^\beta f)\restriction_{\partial\mathbb R^n_+},
\qquad
f\in C_c^\infty(\overline{\mathbb R^n_+}),
\]
admits a unique continuous linear extension
\[
\gamma_0D^\beta\colon
W^{m,p}(\mathbb R^n_+)
\longrightarrow
W^{m-1-|\beta|,p}(\mathbb R^{n-1}).
\]
Moreover,
\[
W^{m,p}_0(\mathbb R^n_+)
=
\left\{
f\in W^{m,p}(\mathbb R^n_+)
\middle|
\gamma_0D^\beta f=0
\text{ for every }|\beta|\leq m-1
\right\},
\]
where \(W^{m,p}_0(\mathbb R^n_+)\) denotes the closure of \(C_c^\infty(\mathbb R^n_+)\) in \(W^{m,p}(\mathbb R^n_+)\).
\end{theorem}

\begin{proof}
Write $x=(x',t)\in\mathbb R^{n-1}\times(0,\infty)$. When $n=1$, interpret $\mathbb R^0$ as a single point equipped with counting measure; the tangential derivatives appearing below then reduce to the derivative of order zero. By Corollary~\ref{teo:extension-sobolev-semiespacio} and density of $C_c^\infty(\mathbb R^n)$ in $W^{m,p}(\mathbb R^n)$, restriction of smooth approximations of an extension of $f$ proves density of $C_c^\infty(\overline{\mathbb R^n_+})$ in $W^{m,p}(\mathbb R^n_+)$.

Begin with $m=1$. If $f\in C_c^\infty(\overline{\mathbb R^n_+})$, the fundamental theorem of calculus gives, for $p=1$,
\[
|f(x',0)|\leq\int_0^\infty|\partial_nf(x',t)|\,dt.
\]
Integrating with respect to $x'$ gives
\[
\|f\restriction_{\partial\mathbb R^n_+}\|_{L^1(\mathbb R^{n-1})}
\leq\|\partial_nf\|_{L^1(\mathbb R^n_+)}.
\]
For $1<p<\infty$, apply the same argument to $|f|^p$ and integrate. Hölder's and Young's inequalities yield
\[
\begin{aligned}
\|f\restriction_{\partial\mathbb R^n_+}\|_{L^p(\mathbb R^{n-1})}^{p}
&\leq p\int_{\mathbb R^n_+}|f|^{p-1}|\partial_nf|\,dx\\
&\leq p\|f\|_{L^p(\mathbb R^n_+)}^{p-1}\|\partial_nf\|_{L^p(\mathbb R^n_+)}\\
&\leq(p-1)\|f\|_{L^p(\mathbb R^n_+)}^p+\|\partial_nf\|_{L^p(\mathbb R^n_+)}^p.
\end{aligned}
\]
By density and completeness of $L^p(\mathbb R^{n-1})$, restriction has a unique continuous linear extension $\gamma_0\colon W^{1,p}(\mathbb R^n_+)\to L^p(\mathbb R^{n-1})$ for every $1\leq p<\infty$.

Now let $m\geq1$, $|\beta|\leq m-1$, and $f\in W^{m,p}(\mathbb R^n_+)$. Each $D^\beta f$ belongs to $W^{1,p}(\mathbb R^n_+)$, so its trace is defined. Take $f_\nu\in C_c^\infty(\overline{\mathbb R^n_+})$ with $f_\nu\to f$ in $W^{m,p}(\mathbb R^n_+)$. For every tangential multi-index $\alpha\in\mathbb N_0^{n-1}$ with $|\alpha|\leq m-1-|\beta|$, the first-order estimate implies
\[
\gamma_0D^{\beta+(\alpha,0)}f_\nu
\longrightarrow
\gamma_0D^{\beta+(\alpha,0)}f
\quad\text{on }L^p(\mathbb R^{n-1}).
\]
For smooth functions,
\[
D_{x'}^\alpha(\gamma_0D^\beta f_\nu)
=\gamma_0D^{\beta+(\alpha,0)}f_\nu.
\]
Integrate this identity against a test function $\varphi\in C_c^\infty(\mathbb R^{n-1})$ and pass to the limit on both sides. Hölder's inequality justifies this passage also for $p=1$, using $\varphi,D^\alpha\varphi\in L^\infty(\mathbb R^{n-1})$. Thus,
\[
D_{x'}^\alpha(\gamma_0D^\beta f)
=\gamma_0D^{\beta+(\alpha,0)}f
\]
in the weak sense. Summing the first-order estimates for these finitely many multi-indices gives
\[
\|\gamma_0D^\beta f\|_{W^{m-1-|\beta|,p}(\mathbb R^{n-1})}
\leq C(n,m,p)\|f\|_{W^{m,p}(\mathbb R^n_+)}.
\]
The density already established also proves uniqueness of this extension.

To characterize the kernel, if $f\in W_0^{m,p}(\mathbb R^n_+)$, approximate by functions in $C_c^\infty(\mathbb R^n_+)$. All their derivatives restrict to zero on the boundary; continuity of the traces implies $\gamma_0D^\beta f=0$ for every $|\beta|\leq m-1$.

Conversely, suppose that these traces vanish. First justify the integration-by-parts identity that allows extension by zero. If $v\in W^{1,p}(\mathbb R^n_+)$ and $\varphi\in C_c^\infty(\mathbb R^n)$, smooth approximation up to the boundary and continuity of $\gamma_0$ give
\[
\int_{\mathbb R^n_+}v\,\partial_i\varphi\,dx
=-\int_{\mathbb R^n_+}(D_iv)\varphi\,dx
-\delta_{in}\int_{\mathbb R^{n-1}}(\gamma_0v)(x')\varphi(x',0)\,dx'.
\]
For smooth functions, the formula follows from the fundamental theorem of calculus and tangential integration by parts; all integrals pass to the limit by Hölder's inequality on the fixed support of $\varphi$. If $\gamma_0v=0$ and $E_0v$ denotes extension by zero to $\mathbb R^n$, the preceding identity is precisely
\[
D_i(E_0v)=E_0(D_iv)
\quad\text{in the sense of distributions on }\mathbb R^n.
\]
Apply this identity successively to $v=D^\alpha f$. At order zero, $E_0f\in L^p(\mathbb R^n)$. If we already know $D^\alpha(E_0f)=E_0(D^\alpha f)$ for some $|\alpha|\leq m-1$, then $D^\alpha f\in W^{1,p}(\mathbb R^n_+)$ and $\gamma_0D^\alpha f=0$. For each $i\in\{1,\dots,n\}$, the first-order step gives
\[
D^{\alpha+e_i}(E_0f)
=D_iE_0(D^\alpha f)
=E_0(D^{\alpha+e_i}f).
\]
Induction proves, for all multi-indices of order at most $m$,
\[
D^\alpha(E_0f)=E_0(D^\alpha f),
\qquad
\|E_0f\|_{W^{m,p}(\mathbb R^n)}=\|f\|_{W^{m,p}(\mathbb R^n_+)}.
\]

Set $F:=E_0f$. For $t>0$, the function $F_t(x):=F(x-te_n)$ belongs to $W^{m,p}(\mathbb R^n)$, vanishes when $x_n<t$, and converges to $F$ in $W^{m,p}(\mathbb R^n)$ as $t\to0^+$. Indeed, $D^\alpha F_t=(D^\alpha F)(\cdot-te_n)$, and continuity of translations in $L^p(\mathbb R^n)$ applies to each of the finitely many derivatives $|\alpha|\leq m$.

For fixed $t$, choose $\chi\in C_c^\infty(\mathbb R^n)$ with $\chi=1$ on a neighborhood of $\overline{B_{\mathrm{euc}}(0,1)}$ and set $\chi_R(x)=\chi(\displaystyle\frac{x}{R})$. The Leibniz rule and $\|D^\alpha\chi_R\|_\infty=R^{-|\alpha|}\|D^\alpha\chi\|_\infty$ show that $\chi_RF_t\to F_t$ in $W^{m,p}(\mathbb R^n)$ as $R\to\infty$: terms with no derivatives on $\chi_R$ converge by domination, while the others contain the factor $R^{-|\alpha|}$, $|\alpha|\geq1$. The product still vanishes on $\{x_n<t\}$. If $\rho_\varepsilon$ is a mollifier supported in $\overline{B_{\mathrm{euc}}(0,\varepsilon)}$ and $0<\varepsilon<\displaystyle\frac{t}{2}$, then
\[
\rho_\varepsilon*(\chi_RF_t)\in C_c^\infty(\mathbb R^n_+),
\]
since its support is contained in $\{x_n\geq t-\varepsilon\}$. For fixed $t,R$, these convolutions converge to $\chi_RF_t$ in $W^{m,p}(\mathbb R^n)$ as $\varepsilon\to0^+$. For each $\nu\in\mathbb N$, choose successively $t_\nu>0$, $R_\nu$, and $0<\varepsilon_\nu<\displaystyle\frac{t_\nu}{2}$ so that each of the three errors is less than $\displaystyle\frac{1}{3\nu}$. The restriction of $\rho_{\varepsilon_\nu}*(\chi_{R_\nu}F_{t_\nu})$ converges to $f$ in $W^{m,p}(\mathbb R^n_+)$. By definition, $f\in W_0^{m,p}(\mathbb R^n_+)$, completing the characterization of the kernel.
\end{proof}

We shall use the following convention from now on. If \(v\in W^{m,p}(B^+)\) satisfies the support assumption in Lemma~\ref{lem:extension-cero-media-bola-semiespacio} and \(\widetilde v\) denotes its extension by zero to \(\mathbb R^n_+\), then for \(|\beta|\leq m-1\) we shall write \(\gamma_0D^\beta v=0\) on \(B'\) to mean that
\[
(\gamma_0D^\beta\widetilde v)\restriction_{B'}=0
\qquad
\text{on }W^{m-1-|\beta|,p}(B').
\]
This notation does not introduce an independent trace on the half-ball; it merely abbreviates the half-space trace after extension by zero.

\begin{lemma}[Locality of the half-space trace]\label{lem:localidad-traza-semiespacio}\index{locality of the half-space trace}
Let $1\leq p<\infty$, let \(O\subset\partial\mathbb R^n_+\) be a relatively open subset, and let \(F\in W^{1,p}(\mathbb R^n_+)\). Suppose that there is an open subset \(V\subset\mathbb R^n\), with \(O\subset V\cap\partial\mathbb R^n_+\), such that \(F=0\) almost everywhere on \(V\cap\mathbb R^n_+\). Then $\gamma_0F=0$ in $O$.
\end{lemma}

\begin{proof}
Fix \(x_0\in O\). Since \(O\) is relatively open in \(\partial\mathbb R^n_+\) and \(V\) is open in \(\mathbb R^n\), we may choose \(r>0\) such that \(B_{\mathrm{euc}}(x_0,2r)\subset V\) and \(B_{\mathrm{euc}}(x_0,2r)\cap\partial\mathbb R^n_+\subset O\). Take a function \(\chi\in C_c^\infty(\mathbb R^n)\) such that \(\operatorname{supp}\chi\subset B_{\mathrm{euc}}(x_0,2r)\) and \(\chi\equiv1\) on \(B_{\mathrm{euc}}(x_0,r)\). Upon restricting \(\chi\) to \(\mathbb R^n_+\), the product \(\chi F\) is well defined as an element of \(W^{1,p}(\mathbb R^n_+)\). Since \(\operatorname{supp}\chi\cap\mathbb R^n_+\subset V\cap\mathbb R^n_+\) and \(F=0\) almost everywhere on \(V\cap\mathbb R^n_+\), we have \(\chi F=0\) in \(W^{1,p}(\mathbb R^n_+)\). Therefore, \(\gamma_0(\chi F)=0\).

We now prove that $\gamma_0(\chi G)=\chi\restriction_{\partial\mathbb R^n_+}\gamma_0G$ for every \(G\in W^{1,p}(\mathbb R^n_+)\). If \(G\in C_c^\infty(\overline{\mathbb R^n_+})\), the identity is immediate, since both sides are the restriction of \(\chi G\) to \(\partial\mathbb R^n_+\). For \(G\in W^{1,p}(\mathbb R^n_+)\), take a sequence \(G_j\in C_c^\infty(\overline{\mathbb R^n_+})\) such that \(G_j\to G\) in \(W^{1,p}(\mathbb R^n_+)\). Since multiplication by the fixed function \(\chi|_{\mathbb R^n_+}\) is continuous on \(W^{1,p}(\mathbb R^n_+)\), we have \(\chi G_j\to\chi G\) in \(W^{1,p}(\mathbb R^n_+)\). Moreover, continuity of the trace gives
\[
\gamma_0(\chi G_j)\longrightarrow \gamma_0(\chi G)
\qquad
\text{on }L^p(\partial\mathbb R^n_+),
\]
and, since multiplication by \(\chi\restriction_{\partial\mathbb R^n_+}\) is continuous on \(L^p(\partial\mathbb R^n_+)\),
\[
\chi\restriction_{\partial\mathbb R^n_+}\gamma_0G_j
\longrightarrow
\chi\restriction_{\partial\mathbb R^n_+}\gamma_0G
\qquad
\text{on }L^p(\partial\mathbb R^n_+).
\]
Passing to the limit in the smooth identity yields
\[
\gamma_0(\chi G)
=
\chi\restriction_{\partial\mathbb R^n_+}\gamma_0G.
\]

Applying this to \(G=F\) gives
\[
0
=
\gamma_0(\chi F)
=
\chi\restriction_{\partial\mathbb R^n_+}\gamma_0F.
\]
Since \(\chi\equiv1\) on \(B_{\mathrm{euc}}(x_0,r)\), it follows that \(\gamma_0F=0\) in \(B_{\mathrm{euc}}(x_0,r)\cap\partial\mathbb R^n_+\). Second countability of $O$ allows us to cover it by a countable family of these neighborhoods. The union of their exceptional sets has measure zero, so $\gamma_0F=0$ almost everywhere on $O$.
\end{proof}

\begin{corollary}[Local characterization on a half-ball]\label{cor:caracterizacion-local-media-bola}\index{local characterization on a half-ball}
Let \(B^+=B_{\mathrm{euc}}(0,R)\cap\mathbb R^n_+\), let \(B'=B_{\mathrm{euc}}(0,R)\cap\partial\mathbb R^n_+\), and let \(v\in W^{m,p}(B^+)\). Suppose that there is a compact set \(K\subset B_{\mathrm{euc}}(0,R)\cap\{x_n\geq0\}\), lying at a positive distance from the spherical portion \(\partial B_{\mathrm{euc}}(0,R)\cap\{x_n\geq0\}\), such that \(v=0\) almost everywhere on \(B^+\setminus K\).

If \(\gamma_0D^\beta v=0\) on \(B'\), in the sense of the preceding convention, for every \(|\beta|\leq m-1\), then \(v\) belongs to the closure of \(C_c^\infty(B^+)\) in \(W^{m,p}(B^+)\).
\end{corollary}

\begin{proof}
If \(K=\varnothing\), then \(v=0\) almost everywhere on \(B^+\), and the result is immediate. Thus assume \(K\neq\varnothing\).

Let \(\widetilde v\) be the extension by zero of \(v\) to \(\mathbb R^n_+\). By Lemma~\ref{lem:extension-cero-media-bola-semiespacio}, we have \(\widetilde v\in W^{m,p}(\mathbb R^n_+)\) and \(D^\beta\widetilde v=\widetilde{D^\beta v}\) for every \(|\beta|\leq m\). By hypothesis, \((\gamma_0D^\beta\widetilde v)\restriction_{B'}=0\) for every \(|\beta|\leq m-1\). We show that in fact \(\gamma_0D^\beta\widetilde v=0\) on all of \(\partial\mathbb R^n_+\).

Fix \(|\beta|\leq m-1\). We already know that the trace vanishes upon restriction to \(B'\). It remains to examine what happens on \(\partial\mathbb R^n_+\setminus B'\). Again denote the spherical portion of the half-ball by \(\Gamma_R:=\partial B_{\mathrm{euc}}(0,R)\cap\{x_n\geq0\}\), take \(\delta:=\operatorname{dist}(K,\Gamma_R)>0\), and set \(U:=B^+\cap\{x\in\mathbb R^n\mid \|x\|>R-\displaystyle\frac{\delta}{2}\}\). As in the proof of the lemma, \(D^\alpha v=0\) almost everywhere on \(U\) for every \(|\alpha|\leq m\).

Let \(x_0=(x_0',0)\in\partial\mathbb R^n_+\setminus B'\). Then \(\|x_0'\|\geq R\). If \(\|x_0'\|>R\), take \(\rho:=\displaystyle\frac{\|x_0'\|-R}{2}>0\). Then \(B_{\mathrm{euc}}(x_0,\rho)\cap\mathbb R^n_+\subset\mathbb R^n_+\setminus B^+\), and the definition of extension by zero gives \(D^\beta\widetilde v=0\) almost everywhere on \(B_{\mathrm{euc}}(x_0,\rho)\cap\mathbb R^n_+\).

If \(\|x_0'\|=R\), take \(\rho:=\displaystyle\frac{\delta}{4}\). If \(x\in B_{\mathrm{euc}}(x_0,\rho)\cap\mathbb R^n_+\), then either \(x\notin B^+\), where \(D^\beta\widetilde v=0\), or \(x\in B^+\). In the latter case, \(\|x\|\geq R-\rho=R-\displaystyle\frac{\delta}{4}>R-\displaystyle\frac{\delta}{2}\), so \(x\in U\), where \(D^\beta v=0\) almost everywhere. Thus in this case too, \(D^\beta\widetilde v=0\) almost everywhere on \(B_{\mathrm{euc}}(x_0,\rho)\cap\mathbb R^n_+\).

Thus, \(D^\beta\widetilde v\) vanishes almost everywhere on a relative neighborhood of each point of \(\partial\mathbb R^n_+\setminus B'\). Since \(|\beta|\leq m-1\), we have \(D^\beta\widetilde v\in W^{1,p}(\mathbb R^n_+)\), and Lemma~\ref{lem:localidad-traza-semiespacio}, applied to \(F=D^\beta\widetilde v\), yields \(\gamma_0D^\beta\widetilde v=0\) in \(\partial\mathbb R^n_+\). Since this holds for every \(|\beta|\leq m-1\), Theorem~\ref{teo:traza-semiespacio} implies \(\widetilde v\in W^{m,p}_0(\mathbb R^n_+)\). Hence there is a sequence \((\eta_k)_{k\in\mathbb N}\subset C_c^\infty(\mathbb R^n_+)\) such that \(\eta_k\to\widetilde v\) in \(W^{m,p}(\mathbb R^n_+)\).

Take a function \(\chi\in C_c^\infty(\mathbb R^n)\) such that \(\operatorname{supp}\chi\subset B_{\mathrm{euc}}(0,R)\) and \(\chi\equiv1\) on a neighborhood of \(K\). Since multiplication by \(\chi|_{\mathbb R^n_+}\) is continuous on \(W^{m,p}(\mathbb R^n_+)\), we have \(\chi\eta_k\to\chi\widetilde v=\widetilde v\) in \(W^{m,p}(\mathbb R^n_+)\). Moreover, \((\chi\eta_k)\restriction_{B^+}\in C_c^\infty(B^+)\), because \(\eta_k\) has compact support in \(\mathbb R^n_+\) and \(\operatorname{supp}\chi\subset B_{\mathrm{euc}}(0,R)\). Restricting to \(B^+\) gives a sequence in \(C_c^\infty(B^+)\) converging to \(v\) in \(W^{m,p}(B^+)\). Therefore, \(v\) belongs to the closure of \(C_c^\infty(B^+)\) in \(W^{m,p}(B^+)\).
\end{proof}

\chapter{Distributions and Fourier analysis}
\label{cap:distribuciones-fourier}

\begin{figura}[htbp]
    \centering
    \includegraphics[height=3.4cm]{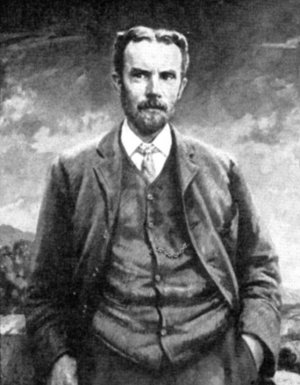}

    \smallskip
    \footnotesize Oliver Heaviside (1850--1925)

    \caption[Oliver Heaviside]{Oliver Heaviside systematically used the step function and operational calculus methods to study electrical problems, anticipating ideas that later received a rigorous formulation in distribution theory.}
    \label{fig:heaviside}
\end{figura}

Distribution theory extends the notion of a function and makes it possible to differentiate objects that have neither a classical derivative nor a weak derivative representable by a locally integrable function. Its modern formulation, due to Laurent Schwartz, rigorously organizes ideas that had already appeared in operational calculus and in the study of discontinuous signals.

Sobolev spaces showed that integration by parts makes it possible to assign derivatives to functions of limited regularity. The distributional framework takes this idea one step further: a distribution is defined by its action on test functions, and its derivatives are obtained by formally transferring differentiation to those functions. This provides a unified treatment of the Heaviside function, the Dirac delta, and many other singularities arising in differential equations, mathematical physics, and Fourier analysis.

A weak solution may therefore have derivatives only in the distributional sense. For differential operators with constant coefficients, the Fourier transform turns those derivatives into multiplication: with the unitary normalization to be fixed below,
\[
\widehat{D^\alpha u}(\xi)
=
(i\xi)^\alpha\widehat u(\xi).
\]
This identity reduces linear equations to algebraic relations in the frequency variable. It also allows us to construct fundamental solutions and regularizations, and to distinguish the behavior of low frequencies from that of high frequencies. Tempered distributions will provide the common domain on which these operations remain well defined.

The step function introduced by Oliver Heaviside provides an immediate motivation. Among other phenomena, it models the instantaneous switching on of a signal and has a jump at the origin; its derivative cannot be represented by a locally integrable function, but it can be represented by a distribution.

\begin{example}[The Heaviside function]\label{ej:distribuciones-y-fourier-la-funcion-de-heaviside}
We define the Heaviside function $H\colon \mathbb{R}\longrightarrow \mathbb{R}$ by
\[
H(x) = \begin{cases}
0, & x < 0,\\
1, & x \geq 0.
\end{cases}
\]
Although $H$ is locally integrable, it is not weakly differentiable in the sense of having its derivative represented by a function in $L^1_{\mathrm{loc}}(\mathbb R)$. This difficulty motivates the distributional extension of the derivative.

Suppose, seeking a contradiction, that there exists $v\in L^{1}_{\mathrm{loc}}(\mathbb{R})$ satisfying
\[
\int_{\mathbb{R}} H(x)\varphi'(x)dx = -\int_{\mathbb{R}} v(x)\varphi(x)dx
\qquad \forall \varphi\in C_c^\infty(\mathbb{R}).
\]
Since $H=0$ on $(-\infty,0)$ and $H=1$ on $[0,\infty)$, we obtain
\[
\int_{\mathbb{R}} H(x)\varphi'(x)dx = \int_0^\infty \varphi'(x)dx.
\]

Since $\varphi$ has compact support, there exists $R>0$ such that $\varphi(x)=0$ for $x\geq R$, so
\[
\int_0^\infty \varphi'(x)dx = \int_0^R \varphi'(x)dx.
\]
This integral is defined over a finite interval, and since $\varphi$ is smooth, the Lebesgue integral agrees with the Riemann integral. Theorem~\ref{teo:b5-fundamental-calculo-riemann-banach} gives
\[
\int_0^R \varphi'(x)dx = \varphi(R)-\varphi(0) = -\varphi(0).
\]

Consequently, the preceding equality would imply
\[
\int_{\mathbb{R}} v(x)\varphi(x)dx = \varphi(0), \qquad \forall \varphi\in C_c^\infty(\mathbb{R}).
\]
Now, if we take a family of functions $\varphi_\varepsilon\in C_c^\infty(\mathbb{R})$ with $0\leq \varphi_\varepsilon \leq 1$, $\varphi_\varepsilon(0)=1$, and $\supp(\varphi_\varepsilon)\subset (-\varepsilon,\varepsilon)$, it follows that
\[
\int_{\mathbb{R}} v(x)\varphi_\varepsilon(x)dx = 1.
\]
But since $v\in L^1_{\text{loc}}(\mathbb{R})$,
\[
\left|\int_{\mathbb{R}} v(x)\varphi_\varepsilon(x)dx\right| \leq \int_{-\varepsilon}^\varepsilon |v(x)|dx \xrightarrow[\varepsilon\to 0]{} 0,
\]
where the limit is justified by applying the \emph{Lebesgue dominated convergence theorem} (Theorem~\ref{convergencia dominada}), since $v\in L^1_{\text{loc}}(\mathbb{R})$ ensures integrability on the intervals $(-\varepsilon,\varepsilon)$. This contradicts the fact that the integral always equals $1$. Thus no such function $v$ exists, and we conclude that $H$ has no weak derivative.
\end{example}

This calculation makes clear the need to extend the notion of a weak derivative: we must allow new objects that, although not ordinary functions, act consistently on test functions. This new object is precisely the \emph{Dirac delta}.

Before formalizing the definition, it is helpful to recall the physical motivation: in his study of quantum mechanics, Paul Dirac (1902–1984) formally introduced the symbol $\delta_{0}$ to model \textit{point sources of energy}. In that setting, the delta was an indispensable tool for describing systems in which all the mass or charge is concentrated at a single point, and engineers and physicists used it long before a rigorous mathematical foundation became available through the distribution theory developed by Laurent Schwartz. We now formalize the Dirac delta.

\begin{example}[The Dirac delta]\label{ej:distribuciones-y-fourier-la-delta-de-dirac}
We define the \textit{\textbf{Dirac delta}} to be the linear functional $\delta_{0}\colon C_{c}^{\infty}(\mathbb{R}^{n})\longrightarrow \mathbb{R}$ given by \[\delta_{0}(\phi)=\phi(0).\]
\end{example}

The preceding identity suggests that $\delta_{0}$ should be interpreted as the derivative of the Heaviside function. It is not an ordinary function, but it can act linearly on test functions. For this action to have analytical stability, linearity alone is insufficient: continuity is needed. Before defining a distribution, we must therefore specify the topology of the space $C_{c}^{\infty}(\mathbb{R}^{n})$.
\section{The space \texorpdfstring{$\mathcal{D}(\Omega)$}{D(Omega)}}\label{sec: espacio D Omega}

Let $\Omega\subseteq\mathbb R^n$ be an open subset. The space $C_c^\infty(\Omega)$ of smooth functions with compact support will be the basic space of test functions. To use it in functional analysis, we need to equip it with a topology that simultaneously records the behavior of a function and all its derivatives, without reducing this information to a single norm.

We first carry out the construction by fixing a compact set $K\subseteq\Omega$ and considering the functions whose support is contained in $K$. We will use the general results on topologies induced by families of seminorms, locally convex spaces, metrizability, and Fréchet spaces collected in Appendix~\ref{ap:evt}.

In particular, we will use Definition~\ref{def: topologia inducida por seminormas}, Theorem~\ref{seminormas topologia}, Propositions~\ref{prop: Hausdorff seminormas separan puntos}, \ref{prop: convergencia y Cauchy por seminormas}, and \ref{prop: metrizacion por seminormas}, as well as Theorem~\ref{teo: caracterizacion Frechet por seminormas}.

We first consider the space of all smooth functions with support contained in a fixed compact subset $K\subseteq\Omega$, which we denote by $\mathcal D_K(\Omega)$:
\[
\mathcal D_K(\Omega):=
\{\phi\in C_c^\infty(\Omega)\mid \operatorname{supp}(\phi)\subseteq K\}.
\]

The topology on this subspace must be consistent with the natural notion of uniform convergence of functions and their partial derivatives of every order on the compact set $K$. For each $j\in\mathbb N_0$, we define
\[
\|\phi\|_{C^j(K)}:=
\max_{|\alpha|\leq j}\sup_{x\in K}|D^\alpha\phi(x)|.
\]
Since the case $\alpha=0$ is included and $\operatorname{supp}(\phi)\subseteq K$, each $\|\cdot\|_{C^j(K)}$ is a norm on $\mathcal D_K(\Omega)$. If $\|\phi\|_{C^j(K)}=0$, then in particular $\displaystyle \|\phi\|_{C^0(K)}=0$, so $\phi=0$ on $K$. On the other hand, since $\operatorname{supp}(\phi)\subseteq K$, if $x\in\Omega\setminus K$, then $\phi(x)=0$. Therefore, $\phi=0$ throughout $\Omega$.

By Theorem~\ref{seminormas topologia}, the family of norms $\{\|\cdot\|_{C^j(K)}\}_{j\in\mathbb N_0}$ induces a topology $\tau_K$ on $\mathcal D_K(\Omega)$ under which this space is a locally convex topological vector space. A neighborhood base at the origin for this topology consists of finite intersections of sets of the form
\[
\bigcap_{q=1}^{k}
\left\{
\phi\in\mathcal D_K(\Omega)
\Biggm|
\|\phi\|_{C^{j_q}(K)}<r_q
\right\},
\]
where $j_1,\dots,j_k\in\mathbb N_0$, $r_1,\dots,r_k>0$, and $k\in\mathbb N$. The neighborhoods of an arbitrary element $\phi\in\mathcal D_K(\Omega)$ have the form $\phi+U$, where $U$ is a neighborhood of the origin.

\begin{lemma}\label{lema: base local numerable DK}
The family
\[
V_{K,j,\ell}:=
\left\{
\phi\in\mathcal D_K(\Omega)
\Biggm|
\|\phi\|_{C^j(K)}<\frac{1}{\ell}
\right\},
\qquad j\in\mathbb N_0,\quad \ell\in\mathbb N,
\]
is a countable neighborhood base at the origin for $\tau_K$. In particular, $(\mathcal D_K(\Omega),\tau_K)$ has a countable local base at each point.
\end{lemma}

\begin{proof}
First, observe that each $V_{K,j,\ell}$ is a neighborhood of the origin, since it is one of the neighborhoods determined by the norm $\|\cdot\|_{C^j(K)}$.

Now let $U$ be a neighborhood of the origin in $\tau_K$. By the description of the local base, there exists a basic neighborhood
\[
W:=
\bigcap_{q=1}^{k}
\left\{
\phi\in\mathcal D_K(\Omega)
\Biggm|
\|\phi\|_{C^{j_q}(K)}<r_q
\right\}
\]
such that $W\subseteq U$. We define $j:=\displaystyle\max\{j_1,\dots,j_k\}$ and choose $\ell\in\mathbb N$ such that $\displaystyle\frac{1}{\ell}<\min\{r_1,\dots,r_k\}$.

Since the norms $\|\cdot\|_{C^j(K)}$ are increasing in $j$, if $\psi\in V_{K,j,\ell}$, then, for each $q\in\{1,\dots,k\}$,
\[
\|\psi\|_{C^{j_q}(K)}
\leq
\|\psi\|_{C^j(K)}
<
\frac{1}{\ell}
<
r_q.
\]
Therefore, $\psi\in W$. Thus, $V_{K,j,\ell}\subseteq W\subseteq U$.

This proves that $\{V_{K,j,\ell}\mid j\in\mathbb N_0,\ \ell\in\mathbb N\}$ is a local base at the origin. Moreover, this family is countable, since $\mathbb N_0\times\mathbb N$ is countable.

Finally, since translations are homeomorphisms in the topology induced by seminorms, for each $\phi\in\mathcal D_K(\Omega)$ the family
\[
\{\phi+V_{K,j,\ell}\mid j\in\mathbb N_0,\ \ell\in\mathbb N\}
\]
is a countable neighborhood base at $\phi$.
\end{proof}

Since the family $\{\|\cdot\|_{C^j(K)}\}_{j\in\mathbb N_0}$ separates points, Proposition~\ref{prop: Hausdorff seminormas separan puntos} implies that the space $(\mathcal D_K(\Omega),\tau_K)$ is Hausdorff. If $\phi\neq0$, then, since $\operatorname{supp}(\phi)\subseteq K$, it is impossible for $\phi$ to vanish throughout $K$; therefore, $\|\phi\|_{C^0(K)}>0$.

In particular, singletons in $(\mathcal D_K(\Omega),\tau_K)$ are closed. Moreover, by Lemma~\ref{lema: base local numerable DK}, this space has a countable local base at the origin. By Theorem~\ref{teo: criterio de metrizabilidad para EVT}, we conclude that $(\mathcal D_K(\Omega),\tau_K)$ is metrizable.

Furthermore, in this case we can write down a compatible metric explicitly. By Proposition~\ref{prop: metrizacion por seminormas}, the topology $\tau_K$ is induced by the translation-invariant metric
\[
d_K(\phi,\psi):=
\displaystyle\sum_{j=0}^{\infty}
\frac{1}{2^j}
\frac{\|\phi-\psi\|_{C^j(K)}}{1+\|\phi-\psi\|_{C^j(K)}},
\qquad
\phi,\psi\in\mathcal D_K(\Omega).
\]

By Proposition~\ref{prop: convergencia y Cauchy por seminormas}, a sequence $(\phi_n)_{n\in\mathbb N}$ converges to $\phi$ in $\mathcal D_K(\Omega)$ if and only if
\[
\|\phi_n-\phi\|_{C^j(K)}\to0
\qquad
\text{for every }j\in\mathbb N_0.
\]
This means that $\phi_n\to\phi$ if and only if, for every multi-index $\alpha$, the sequence $(D^\alpha\phi_n)_{n\in\mathbb N}$ converges uniformly on $K$ to $D^\alpha\phi$. If $\|\phi_n-\phi\|_{C^j(K)}\to0$ for all $j$, then, given a multi-index $\alpha$, it suffices to take $j\geq |\alpha|$ to obtain
\[
\|D^\alpha\phi_n-D^\alpha\phi\|_{C^0(K)}
\leq
\|\phi_n-\phi\|_{C^j(K)}\to0.
\]
Conversely, if for every multi-index $\alpha$ the sequence $(D^\alpha\phi_n)_{n\in\mathbb N}$ converges uniformly on $K$ to $D^\alpha\phi$, then, for each $j\in\mathbb N_0$, the set of multi-indices $\alpha$ with $|\alpha|\leq j$ is finite. Therefore,
\[
\|\phi_n-\phi\|_{C^j(K)}
=
\max_{|\alpha|\leq j}
\|D^\alpha\phi_n-D^\alpha\phi\|_{C^0(K)}
\to0.
\]

Likewise, a sequence $(\phi_n)_{n\in\mathbb N}$ is Cauchy in $\mathcal D_K(\Omega)$ if and only if
\[
\|\phi_n-\phi_m\|_{C^j(K)}\to0
\qquad
\text{as }n,m\to\infty,
\]
for every $j\in\mathbb N_0$. Equivalently, for every multi-index $\alpha$, the sequence $(D^\alpha\phi_n)_{n\in\mathbb N}$ is uniformly Cauchy on $K$.

\section{Projective limits}\label{sec: limites proyectivos}

A natural question is whether the space $\mathcal D_K(\Omega)$ is complete under the topology $\tau_K$ we have just constructed. It is indeed complete. Moreover, $\mathcal D_K(\Omega)$ is not only a locally convex space generated by a countable family of norms, but can also be described as a \textit{projective limit of Banach spaces}.

The same description by projective limits will be used again for spaces of test sections of smooth vector bundles and for the space of rapidly decreasing functions, usually called the \textit{Schwartz class}. The latter will be the test function space for tempered distributions.

From the topological point of view, the projective topology we will use is not a new construction: it is precisely an initial topology in the sense of Definition~\ref{def: topologia inicial}. We call it \textit{projective} because it is associated with a family of linear maps into known topological spaces.

\begin{definition}[Projective topology]
\label{def: topologia proyectiva}\index{projective topology}
Let $X$ be a vector space, let $\{X_i\}_{i\in I}$ be a family of topological vector spaces, and, for each $i\in I$, let $T_i\colon X\longrightarrow X_i$ be a linear map. The \textit{projective topology} on $X$ associated with the family $(T_i)_{i\in I}$ is the initial topology determined by these maps.

Equivalently, it is the coarsest topology on $X$ for which all the maps $T_i\colon X\longrightarrow X_i$ are continuous. A neighborhood base at the origin for this topology consists of the sets
\[
\bigcap_{q=1}^{r}T_{i_q}^{-1}(U_{i_q}),
\]
where $r\in\mathbb N$, $i_1,\dots,i_r\in I$, and each $U_{i_q}$ is a neighborhood of the origin in $X_{i_q}$.
\end{definition}

\begin{proposition}
\label{prop: propiedades topologia proyectiva}
Let $X$ be a vector space, let $\{X_i\}_{i\in I}$ be a family of topological vector spaces, and, for each $i\in I$, let $T_i\colon X\longrightarrow X_i$ be a linear map. Then the projective topology associated with $(T_i)_{i\in I}$ makes $X$ a topological vector space.

Moreover:
\begin{enumerate}[label=(\alph*)]
 \item if each $X_i$ is locally convex, then $X$ is locally convex;

 \item if each $X_i$ is Hausdorff and the family $(T_i)_{i\in I}$ separates points of $X$, then $X$ is Hausdorff.
\end{enumerate}
\end{proposition}

\begin{proof}
Denote the projective topology on $X$ by $\tau_{\mathrm{proj}}$.

We first prove that addition is continuous. Since it suffices to check continuity on basic neighborhoods, let $x,y\in X$ and let $W$ be a basic neighborhood of $x+y$ in $\tau_{\mathrm{proj}}$. We can write
\[
W=
\bigcap_{q=1}^{r}T_{i_q}^{-1}(U_q),
\]
where, for each $q\in\{1,\dots,r\}$, the set $U_q$ is a neighborhood of $T_{i_q}(x+y)$ in $X_{i_q}$.

Since addition is continuous in $X_{i_q}$, there exist neighborhoods $V_q$ of $T_{i_q}(x)$ and $W_q$ of $T_{i_q}(y)$ such that $V_q+W_q\subseteq U_q$. We define
\[
V:=
\bigcap_{q=1}^{r}T_{i_q}^{-1}(V_q),
\qquad
W':=
\bigcap_{q=1}^{r}T_{i_q}^{-1}(W_q).
\]
Then $V$ is a neighborhood of $x$ and $W'$ is a neighborhood of $y$. If $u\in V$ and $v\in W'$, then, for each $q\in\{1,\dots,r\}$, we have $T_{i_q}(u)\in V_q$ and $T_{i_q}(v)\in W_q$. Since $T_{i_q}$ is linear,
\[
T_{i_q}(u+v)=T_{i_q}(u)+T_{i_q}(v)\in V_q+W_q\subseteq U_q.
\]
Therefore, $u+v\in W$. Thus, $V+W'\subseteq W$, and addition is continuous.

We now show that scalar multiplication is continuous. Since it suffices to check continuity on basic neighborhoods, let $\lambda\in\mathbb K$, $x\in X$, and let $W$ be a basic neighborhood of $\lambda x$. We write
\[
W=
\bigcap_{q=1}^{r}T_{i_q}^{-1}(U_q),
\]
where $U_q$ is a neighborhood of $\lambda T_{i_q}(x)$ in $X_{i_q}$. Since scalar multiplication is continuous in $X_{i_q}$, there exist a neighborhood $A_q$ of $\lambda$ in $\mathbb K$ and a neighborhood $V_q$ of $T_{i_q}(x)$ in $X_{i_q}$ such that $A_qV_q\subseteq U_q$. We define
\[
A:=\bigcap_{q=1}^{r}A_q,
\qquad
V:=\bigcap_{q=1}^{r}T_{i_q}^{-1}(V_q).
\]
Then $A$ is a neighborhood of $\lambda$ in $\mathbb K$ and $V$ is a neighborhood of $x$ in $X$. If $\mu\in A$ and $u\in V$, then, for each $q\in\{1,\dots,r\}$, we have $\mu\in A_q$ and $T_{i_q}(u)\in V_q$. Therefore,
\[
T_{i_q}(\mu u)=\mu T_{i_q}(u)\in A_qV_q\subseteq U_q.
\]
Thus, $\mu u\in W$. We conclude that scalar multiplication is continuous.

This proves that $(X,\tau_{\mathrm{proj}})$ is a topological vector space.

We prove (a). Suppose that each $X_i$ is locally convex. For each $i\in I$, take a neighborhood base at the origin in $X_i$ consisting of convex balanced sets. Then the basic neighborhoods of the origin in $X$ may be taken to have the form
\[
\bigcap_{q=1}^{r}T_{i_q}^{-1}(U_{i_q}),
\]
where each $U_{i_q}$ is convex and balanced. Since the inverse image of a convex balanced set under a linear map retains both properties, and since finite intersections of convex balanced sets are also convex and balanced, we conclude that $X$ is locally convex.

We prove (b). Suppose that each $X_i$ is Hausdorff and that the family $(T_i)_{i\in I}$ separates points of $X$. Then
\[
\bigcap_{i\in I}\ker(T_i)=\{0\}.
\]
Since $X_i$ is Hausdorff, the set $\{0\}$ is closed in $X_i$. Since $T_i$ is continuous, $\ker(T_i)=T_i^{-1}(\{0\})$ is closed in $X$. Therefore, $\{0\}$ is closed in $X$, being an intersection of closed sets. By Proposition~\ref{prop: Hausdorff si singleton cerrado}, we conclude that $X$ is Hausdorff.
\end{proof}

The general formulation of projective limits uses systems of spaces and transition maps. In the applications in this book, we will mainly encounter the concrete case in which a space is realized as an intersection of Banach spaces contained in a common topological vector space.

\begin{definition}[Projective limit realized as an intersection]
\label{def: limite proyectivo como interseccion}\index{projective limit realized as an intersection}
Let $Y$ be a Hausdorff topological vector space and let $\{X_n\}_{n\in\mathbb N}$ be a family of Banach spaces such that each $X_n$ is a vector subspace of $Y$ and the inclusion $X_n\hookrightarrow Y$ is continuous. We define
\[
X:=
\bigcap_{n\in\mathbb N}X_n.
\]
For each $n\in\mathbb N$, denote the natural inclusion by $\iota_n\colon X\longrightarrow X_n$ and equip $X$ with the projective topology associated with the family $(\iota_n)_{n\in\mathbb N}$.

In this setting we write
\[
X=
\varprojlim_{n\in\mathbb N}X_n
\]
and say that $X$ is the \textit{projective limit} of the family $\{X_n\}_{n\in\mathbb N}$. Some texts also use the term \textit{inverse limit}.
\end{definition}

\begin{proposition}[Sequential completeness of a projective intersection]
\label{prop: completitud secuencial limite proyectivo}\index{sequential completeness of a projective intersection}
Let $Y$ be a Hausdorff topological vector space and let $\{X_n\}_{n\in\mathbb N}$ be a family of Banach spaces such that each $X_n$ is a vector subspace of $Y$ and the inclusion $X_n\hookrightarrow Y$ is continuous. Consider
\[
X=
\bigcap_{n\in\mathbb N}X_n
\]
equipped with the projective topology associated with the inclusions $\iota_n\colon X\longrightarrow X_n$. Then $X$ is sequentially complete: every Cauchy sequence in $X$ converges in $X$.
\end{proposition}

\begin{proof}
Let $(x^{(q)})_{q\in\mathbb N}$ be a Cauchy sequence in $X$ with respect to the projective topology. Fix $n\in\mathbb N$. We will prove that $(x^{(q)})_{q\in\mathbb N}$ is a Cauchy sequence in $X_n$.

Let $U_n$ be a neighborhood of the origin in $X_n$. Since $\iota_n\colon X\longrightarrow X_n$ is continuous by definition of the projective topology, the set $\iota_n^{-1}(U_n)$ is a neighborhood of the origin in $X$. Since $(x^{(q)})_{q\in\mathbb N}$ is Cauchy in $X$, there exists $q_0\in\mathbb N$ such that
\[
x^{(q)}-x^{(p)}\in\iota_n^{-1}(U_n)
\qquad
\forall q,p\geq q_0.
\]
Therefore,
\[
x^{(q)}-x^{(p)}\in U_n
\qquad
\forall q,p\geq q_0.
\]
This proves that $(x^{(q)})_{q\in\mathbb N}$ is Cauchy in $X_n$.

Since $X_n$ is Banach, there exists $x_n\in X_n$ such that $x^{(q)}\to x_n$ in $X_n$ as $q\to\infty$.

We show that the limits $x_n$ do not depend on $n$ when regarded as elements of $Y$. Let $n,m\in\mathbb N$. Since the inclusions $X_n\hookrightarrow Y$ and $X_m\hookrightarrow Y$ are continuous, $x^{(q)}\to x_n$ in $X_n$ implies $x^{(q)}\to x_n$ in $Y$, and $x^{(q)}\to x_m$ in $X_m$ implies $x^{(q)}\to x_m$ in $Y$. Since $Y$ is Hausdorff, the limit of a convergent sequence is unique. Consequently, $x_n=x_m$ in $Y$.

There is therefore a unique element $x\in Y$ such that $x=x_n$ for every $n\in\mathbb N$. Since $x_n\in X_n$ for each $n$, we conclude that $x\in\displaystyle\bigcap_{n\in\mathbb N}X_n=X$.

It remains to prove that $x^{(q)}\to x$ in the projective topology of $X$. Let $V$ be a neighborhood of the origin in $X$. By definition of this topology, there is a basic neighborhood
\[
V_0=
\bigcap_{s=1}^{r}\iota_{n_s}^{-1}(U_{n_s})
\]
such that $V_0\subseteq V$, where $r\in\mathbb N$, $n_1,\dots,n_r\in\mathbb N$, and each $U_{n_s}$ is a neighborhood of the origin in $X_{n_s}$.

For each $s\in\{1,\dots,r\}$, we know that $x^{(q)}\to x$ in $X_{n_s}$. Therefore, there exists $q_s\in\mathbb N$ such that $x^{(q)}-x\in U_{n_s}$ for every $q\geq q_s$. Taking $q_*:=\displaystyle\max\{q_1,\dots,q_r\}$, we obtain, for every $q\geq q_*$,
\[
x^{(q)}-x\in
\bigcap_{s=1}^{r}\iota_{n_s}^{-1}(U_{n_s})
=
V_0
\subseteq V.
\]
Consequently, $x^{(q)}\to x$ in $X$. We conclude that $X$ is sequentially complete.
\end{proof}

\begin{remark}\label{obs:distribuciones-y-fourier-completitud-espacio-vectorial-topologico-formula-mediante}
If completeness of a topological vector space is formulated in terms of Cauchy nets or filters, the same argument shows that $\displaystyle\varprojlim_{n\in\mathbb N}X_n$ is complete in that sense. In the immediate applications, the sequential formulation will suffice, since the spaces we will study are metrizable.
\end{remark}

We now return to the space $\mathcal D_K(\Omega)$. For each $j\in\mathbb N_0$, we define
\[
C_K^j(\Omega):=
\left\{
\phi\in C^j(\Omega)
\ \middle|\
\operatorname{supp}(\phi)\subseteq K
\right\},
\]
and equip it with the norm
\[
\|\phi\|_{C^j(K)}
:=
\max_{|\alpha|\leq j}\sup_{x\in K}|D^\alpha\phi(x)|.
\]

\begin{lemma}
\label{lema: CKj es Banach}
For each $j\in\mathbb N_0$, the space $\bigl(C_K^j(\Omega),\|\cdot\|_{C^j(K)}\bigr)$ is a Banach space.
\end{lemma}

\begin{proof}
Let $(\phi_n)_{n\in\mathbb N}$ be a Cauchy sequence in $C_K^j(\Omega)$. Fix a multi-index $\alpha$ with $|\alpha|\leq j$.

If $x\in\Omega\setminus K$, then, since $\operatorname{supp}(\phi_n)\subseteq K$, we have $x\notin\operatorname{supp}(\phi_n)$. Therefore, there is an open neighborhood of $x$ on which $\phi_n$ vanishes. Consequently, $D^\alpha\phi_n(x)=0$. Thus, $D^\alpha\phi_n$ vanishes on $\Omega\setminus K$ for every $n\in\mathbb N$.

Moreover, for any $n,m\in\mathbb N$,
\[
\|D^\alpha\phi_n-D^\alpha\phi_m\|_{L^\infty(\Omega)}
=
\|D^\alpha\phi_n-D^\alpha\phi_m\|_{C^0(K)}
\leq
\|\phi_n-\phi_m\|_{C^j(K)}.
\]
Since $(\phi_n)_{n\in\mathbb N}$ is Cauchy in $C_K^j(\Omega)$, it follows that $(D^\alpha\phi_n)_{n\in\mathbb N}$ is uniformly Cauchy on $\Omega$.

We now define the corresponding limit. For each $x\in\Omega$, the real sequence $(D^\alpha\phi_n(x))_{n\in\mathbb N}$ is Cauchy and therefore converges in $\mathbb R$. We define
\[
g_\alpha(x):=
\lim_{n\to\infty}D^\alpha\phi_n(x),
\qquad x\in\Omega.
\]
We show that the convergence is uniform. Let $\varepsilon>0$. Since $(D^\alpha\phi_n)_{n\in\mathbb N}$ is uniformly Cauchy, there exists $N\in\mathbb N$ such that, if $n,m\geq N$, then
\[
\|D^\alpha\phi_n-D^\alpha\phi_m\|_{L^\infty(\Omega)}
<
\varepsilon.
\]
Fix $n\geq N$ and let $m\to\infty$. For each $x\in\Omega$,
\[
|D^\alpha\phi_n(x)-g_\alpha(x)|
\leq
\varepsilon.
\]
Taking the supremum over $x\in\Omega$, we obtain
\[
\|D^\alpha\phi_n-g_\alpha\|_{L^\infty(\Omega)}
\leq
\varepsilon.
\]
Therefore, $D^\alpha\phi_n\to g_\alpha$ uniformly on $\Omega$.

Since each $D^\alpha\phi_n$ is continuous and the convergence is uniform, $g_\alpha$ is continuous on $\Omega$. Moreover, since $D^\alpha\phi_n=0$ on $\Omega\setminus K$ for every $n$, passing to the limit gives
\[
g_\alpha=0
\qquad
\text{in }\Omega\setminus K.
\]

We define $\phi:=g_0$. We will prove that $\phi\in C_K^j(\Omega)$ and that $D^\alpha\phi=g_\alpha$ for every multi-index $\alpha$ with $|\alpha|\leq j$.

We begin by showing compatibility between the limits. We claim that, if $\beta$ is a multi-index with $|\beta|\leq j-1$ and $i\in\{1,\dots,n\}$, then
\[
\partial_i g_\beta=g_{\beta+e_i}
\qquad
\text{in }\Omega.
\]

Fix $x\in\Omega$. Since $\Omega$ is open, there exists $\delta>0$ such that $x+te_i\in\Omega$ for every $t\in[-\delta,\delta]$. If $h\in(-\delta,\delta)$, then, for each $n\in\mathbb N$, Theorem~\ref{teo:b5-fundamental-calculo-riemann-banach} applied to the function $t\mapsto D^\beta\phi_n(x+te_i)$ gives
\[
D^\beta\phi_n(x+he_i)-D^\beta\phi_n(x)
=
\int_0^h D^{\beta+e_i}\phi_n(x+te_i)\,dt.
\]

Let $n\to\infty$. By the uniform convergence $D^\beta\phi_n\to g_\beta$, the left-hand side converges to
\[
g_\beta(x+he_i)-g_\beta(x).
\]
On the other hand,
\[
\left|
\int_0^h
\bigl(
D^{\beta+e_i}\phi_n(x+te_i)
-
g_{\beta+e_i}(x+te_i)
\bigr)
\,dt
\right|
\leq
|h|
\|D^{\beta+e_i}\phi_n-g_{\beta+e_i}\|_{L^\infty(\Omega)}.
\]
Since $D^{\beta+e_i}\phi_n\to g_{\beta+e_i}$ uniformly on $\Omega$, the right-hand side converges to $0$. Therefore,
\[
g_\beta(x+he_i)-g_\beta(x)
=
\int_0^h g_{\beta+e_i}(x+te_i)\,dt.
\]

If $h\neq0$, we divide by $h$:
\[
\frac{g_\beta(x+he_i)-g_\beta(x)}{h}
=
\frac{1}{h}
\int_0^h g_{\beta+e_i}(x+te_i)\,dt.
\]
Since $g_{\beta+e_i}$ is continuous at $x$, we have
\[
\left|
\frac{1}{h}
\int_0^h g_{\beta+e_i}(x+te_i)\,dt
-
g_{\beta+e_i}(x)
\right|
\leq
\sup_{t\in[-|h|,|h|]}
|g_{\beta+e_i}(x+te_i)-g_{\beta+e_i}(x)|
\longrightarrow 0
\]
as $h\to0$. Consequently,
\[
\partial_i g_\beta(x)=g_{\beta+e_i}(x).
\]
Since $x\in\Omega$ was arbitrary, we conclude that $\partial_i g_\beta=g_{\beta+e_i}$ on $\Omega$.

Using this relation, we will prove by induction on $r\in\{0,\dots,j\}$ that:
\[
\phi\in C^r(\Omega)
\qquad\text{and}\qquad
D^\alpha\phi=g_\alpha
\quad
\text{for every }|\alpha|\leq r.
\]

For $r=0$, this is immediate, since $\phi=g_0$ and $g_0$ is continuous.

Now suppose that the assertion holds for some $r\in\{0,\dots,j-1\}$. We will prove that it holds for $r+1$. Let $\alpha$ be a multi-index with $|\alpha|=r+1$. Choose $i\in\{1,\dots,n\}$ such that $\alpha_i\geq1$ and write $\alpha=\beta+e_i$, where $|\beta|=r$.

By the induction hypothesis, $D^\beta\phi=g_\beta$. Since $|\beta|=r\leq j-1$, the identity proved above gives
\[
\partial_i g_\beta=g_{\beta+e_i}=g_\alpha.
\]
Therefore,
\[
D^\alpha\phi
=
\partial_i(D^\beta\phi)
=
\partial_i g_\beta
=
g_\alpha.
\]
Moreover, $g_\alpha$ is continuous. Since this holds for every multi-index $\alpha$ with $|\alpha|=r+1$, we conclude that $\phi\in C^{r+1}(\Omega)$ and that $D^\alpha\phi=g_\alpha$ for every $|\alpha|\leq r+1$.

By induction, $\phi\in C^j(\Omega)$ and $D^\alpha\phi=g_\alpha$ for every multi-index $\alpha$ with $|\alpha|\leq j$.

Since $g_0=0$ on $\Omega\setminus K$, we have $\phi=0$ on $\Omega\setminus K$. Since $K$ is closed in $\Omega$, we conclude that $\operatorname{supp}(\phi)\subseteq K$. Therefore, $\phi\in C_K^j(\Omega)$.

Finally, for every multi-index $\alpha$ with $|\alpha|\leq j$, we already know that $D^\alpha\phi_n\to D^\alpha\phi$ uniformly on $\Omega$, and in particular uniformly on $K$. Since the set of multi-indices $\alpha$ with $|\alpha|\leq j$ is finite,
\[
\|\phi_n-\phi\|_{C^j(K)}
=
\max_{|\alpha|\leq j}
\|D^\alpha\phi_n-D^\alpha\phi\|_{C^0(K)}
\longrightarrow 0.
\]
Thus every Cauchy sequence in $C_K^j(\Omega)$ converges in $C_K^j(\Omega)$. We conclude that $\bigl(C_K^j(\Omega),\|\cdot\|_{C^j(K)}\bigr)$ is a Banach space.
\end{proof}

Now observe that
\[
\mathcal D_K(\Omega)
=
\bigcap_{j\in\mathbb N_0}C_K^j(\Omega).
\]
If $\phi\in\mathcal D_K(\Omega)$, then $\phi\in C_K^j(\Omega)$ for every $j\in\mathbb N_0$. Conversely, if $\phi\in\displaystyle\bigcap_{j\in\mathbb N_0}C_K^j(\Omega)$, then $\phi$ is of class $C^j$ for every $j\in\mathbb N_0$, so $\phi\in C^\infty(\Omega)$. Moreover, $\operatorname{supp}(\phi)\subseteq K$. Therefore, $\phi\in\mathcal D_K(\Omega)$.

For each $j\in\mathbb N_0$, consider the natural inclusion $\iota_j\colon \mathcal D_K(\Omega)\longrightarrow C_K^j(\Omega)$. Equip $\mathcal D_K(\Omega)$ with the projective topology associated with the family $(\iota_j)_{j\in\mathbb N_0}$, and denote it by $\tau_K^{\mathrm{proj}}$.

A neighborhood base at the origin for $\tau_K^{\mathrm{proj}}$ consists of finite intersections of sets of the form
\[
\iota_j^{-1}
\left(
\left\{
\phi\in C_K^j(\Omega)
\ \middle|\
\|\phi\|_{C^j(K)}<r
\right\}
\right),
\qquad
j\in\mathbb N_0,\quad r>0.
\]
Since $\iota_j$ is the natural inclusion, these sets agree with
\[
\left\{
\phi\in\mathcal D_K(\Omega)
\ \middle|\
\|\phi\|_{C^j(K)}<r
\right\}.
\]
Therefore, a neighborhood base at the origin for $\tau_K^{\mathrm{proj}}$ consists of the sets
\[
\bigcap_{q=1}^{k}
\left\{
\phi\in\mathcal D_K(\Omega)
\ \middle|\
\|\phi\|_{C^{j_q}(K)}<r_q
\right\},
\]
where $k\in\mathbb N$, $j_1,\dots,j_k\in\mathbb N_0$, and $r_1,\dots,r_k>0$.

This family is the neighborhood base at the origin defining the topology $\tau_K$ constructed earlier from the family of norms $\{\|\cdot\|_{C^j(K)}\}_{j\in\mathbb N_0}$. Consequently, $\tau_K^{\mathrm{proj}}=\tau_K$.

In other words,
\[
\bigl(\mathcal D_K(\Omega),\tau_K\bigr)
=
\varprojlim_{j\in\mathbb N_0}
C_K^j(\Omega).
\]

Since each $C_K^j(\Omega)$ is a Banach space and the family is indexed by $\mathbb N_0$, we can reindex it using $n\mapsto n-1$ to write it as a family indexed by $\mathbb N$. Proposition~\ref{prop: completitud secuencial limite proyectivo} then implies that $\mathcal D_K(\Omega)$ is sequentially complete with respect to $\tau_K$.

It remains to translate this sequential completeness into completeness of the metric $d_K$ already constructed. Let $(\phi_n)_{n\in\mathbb N}$ be a Cauchy sequence with respect to $d_K$. Fix $j\in\mathbb N_0$ and $\varepsilon>0$. Since $\displaystyle\frac{1}{2^j}\displaystyle\frac{\varepsilon}{1+\varepsilon}>0$, there exists $N\in\mathbb N$ such that, if $n,m\geq N$, then
\[
d_K(\phi_n,\phi_m)
<
\frac{1}{2^j}
\frac{\varepsilon}{1+\varepsilon}.
\]
By definition of $d_K$,
\[
\frac{1}{2^j}
\frac{\|\phi_n-\phi_m\|_{C^j(K)}}{1+\|\phi_n-\phi_m\|_{C^j(K)}}
\leq
d_K(\phi_n,\phi_m).
\]
Therefore,
\[
\frac{\|\phi_n-\phi_m\|_{C^j(K)}}{1+\|\phi_n-\phi_m\|_{C^j(K)}}
<
\frac{\varepsilon}{1+\varepsilon}.
\]
Since the function $t\mapsto \displaystyle\frac{t}{1+t}$ is strictly increasing on $[0,\infty)$, we conclude that $\|\phi_n-\phi_m\|_{C^j(K)}<\varepsilon$ for any $n,m\geq N$.

We have proved that, for each $j\in\mathbb N_0$, $\|\phi_n-\phi_m\|_{C^j(K)}\to0$ as $n,m\to\infty$. By Proposition~\ref{prop: convergencia y Cauchy por seminormas}, $(\phi_n)_{n\in\mathbb N}$ is Cauchy with respect to $\tau_K$. Since $\mathcal D_K(\Omega)$ is sequentially complete for this topology, there exists $\phi\in\mathcal D_K(\Omega)$ such that $\phi_n\to\phi$ in $\tau_K$.

Finally, since the metric $d_K$ induces the topology $\tau_K$, it follows that $d_K(\phi_n,\phi)\to0$. Thus every Cauchy sequence with respect to $d_K$ converges in $\mathcal D_K(\Omega)$, and therefore $d_K$ is complete.

We conclude that $(\mathcal D_K(\Omega),\tau_K)$ is a locally convex, Hausdorff, metrizable, and complete space. By Theorem~\ref{teo: caracterizacion Frechet por seminormas}, $\mathcal D_K(\Omega)$ is a Fréchet space.

\section{The topology of $\mathcal{D}(\Omega)$}
The construction in this section follows the treatment in \cite[Section~10.1]{Leoni2017}. To construct a topology on $C_{c}^{\infty}(\Omega)$, let $\mathcal{B}_{0}$ be the collection of subsets $U\subseteq C_{c}^{\infty}(\Omega)$ such that:

\begin{enumerate}
\item $0\in U$,
\item $U$ is balanced\footnote{A subset $U$ of a vector space is said to be balanced if for every scalar $\lambda$ with $|\lambda|\leq 1$ we have $\lambda U\subseteq U$.} and convex,
\item for every compact set $K\subseteq \Omega$, the set $U\cap \mathcal{D}_{K}(\Omega)$ is a neighborhood of the origin in $(\mathcal{D}_{K}(\Omega),\tau_{K})$.
\end{enumerate}

The collection $\mathcal{B}_{0}$ is designed to serve as a neighborhood base at the origin in $C_{c}^{\infty}(\Omega)$. From it, we define
\[
\mathcal{B}:=\{\phi+V\mid \phi\in C_{c}^{\infty}(\Omega),\ V\in \mathcal{B}_{0}\}.
\]

These conditions determine the topology through its neighborhoods.

\begin{theorem}\label{teo: topologia de D(Omega)}\index{test function space!topology of \(\mathcal D(\Omega)\)}
Let $\Omega\subseteq \mathbb{R}^{n}$ be open. There is a unique topology $\tau_{LF}$ on $C_c^\infty(\Omega)$ for which $\mathcal B_0$ is a neighborhood base at the origin and the vector space operations are continuous. This topology is Hausdorff and locally convex. For each $\phi\in C_c^\infty(\Omega)$, the family $\phi+\mathcal B_0$ is a neighborhood base at $\phi$; a set $O$ is open if and only if for each $\phi\in O$ there exists $V\in\mathcal B_0$ with $\phi+V\subseteq O$. We denote the topological vector space $(C_c^\infty(\Omega),\tau_{LF})$ by $\mathcal D(\Omega)$ and call its elements \textit{test functions}.
\end{theorem}

\begin{proof}
The family $\mathcal B_0$ is nonempty, since it contains the entire space. Each $U\in\mathcal B_0$ is absorbing: if $\phi\in C_c^\infty(\Omega)$, take a compact set $K$ containing its support. The neighborhood $U\cap\mathcal D_K(\Omega)$ absorbs $\phi$ in that space, and hence $U$ also absorbs it in $C_c^\infty(\Omega)$. Moreover, $\mathcal B_0$ is stable under finite intersections and multiplication by positive scalars. Indeed, these operations preserve convexity and balancedness, and the intersections with each $\mathcal D_K(\Omega)$ remain neighborhoods of the origin. Proposition~\ref{prop: base-discos-topologia-lc} therefore yields a unique locally convex topology with neighborhood base $\mathcal B_0$.

To check that it is Hausdorff, consider
\[
V_\varepsilon:=\left\{\phi\in C_c^\infty(\Omega)
\ \middle|\ \|\phi\|_{L^\infty(\Omega)}<\varepsilon\right\},
\qquad \varepsilon>0.
\]
This set belongs to $\mathcal B_0$, because its intersection with $\mathcal D_K(\Omega)$ is the open ball of radius $\varepsilon$ for the uniform norm. If $\phi\neq\psi$, choose $\displaystyle 3\varepsilon<\|\phi-\psi\|_{L^\infty(\Omega)}$. The neighborhoods $\phi+V_\varepsilon$ and $\psi+V_\varepsilon$ are disjoint: otherwise, the triangle inequality would give a uniform distance less than $2\varepsilon$ between $\phi$ and $\psi$. Their interiors are disjoint open neighborhoods. The description of the open sets and uniqueness follow from Proposition~\ref{prop: topologia-determinada-por-base-local}.
\end{proof}

The norms introduced earlier depend on the compact set $K$ containing the support of the functions. A natural question is whether it is possible to avoid this dependence and define a topology directly on $C_{c}^{\infty}(\Omega)$ using global norms.

With this in mind, we could consider, for each $j\in \mathbb{N}_{0}$, the function
\begin{equation}\label{eq: normas globales en D(Omega)}
\|\phi\|_{j}:=\sup_{\substack{x\in \Omega\\|\alpha|\leq j}}|D^{\alpha}\phi(x)|.
\end{equation}
Observe that this supremum is finite, since every function $\phi\in C_{c}^{\infty}(\Omega)$ has compact support. Moreover, each $\|\cdot\|_{j}$ defines a norm on $C_{c}^{\infty}(\Omega)$.

We could then try to define a topology on $C_{c}^{\infty}(\Omega)$ as the topology induced by the family of norms $\{\|\cdot \|_{j}\}_{j\in \mathbb{N}_{0}}$. However, this construction is unsuitable. Intuitively, these norms do not uniformly control the behavior of the supports of the functions, which leads to a topology that is too weak.

In particular, the resulting space is not complete, as the following example shows.
\begin{example}\label{ej:distribuciones-y-fourier-funcion-tal-sucesion-funciones-funcion-montana}
Consider $\Omega = \mathbb{R}$, and let $\phi \in C_{c}^{\infty}(\mathbb{R})$ be a function such that $\operatorname{supp}(\phi) \subseteq [0,1]$ and $\phi(x) > 0$ for every $x \in (0,1)$. Define the sequence of functions $(\phi_{k})_{k\in \mathbb{N}}$ by
\[
\phi_{k}(x) := \displaystyle\sum_{m=1}^{k} \frac{1}{m}\phi(x-m)
= \phi(x-1) + \frac{1}{2}\phi(x-2) + \dots + \frac{1}{k}\phi(x-k).
\]
\begin{center}
 \includegraphics[scale=0.75]{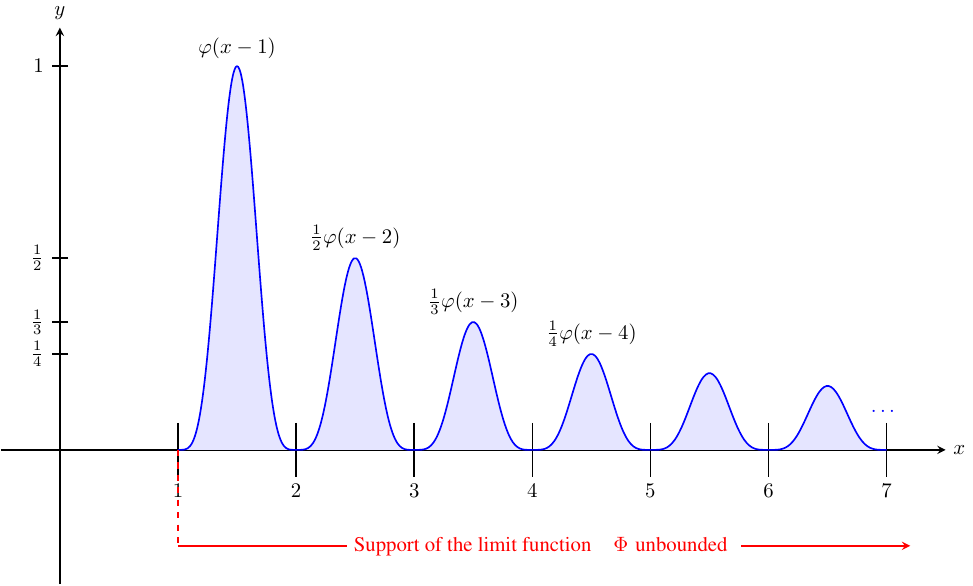}
\end{center}
We will show that this is a Cauchy sequence with respect to the topology generated by the family of seminorms
\[
\|\cdot\|_{j} = \sup_{\substack{x \in \mathbb{R}\\ i \leq j}} |D^{i}(\cdot)(x)|,
\]
but its limit does not have compact support, concluding that $C_{c}^{\infty}(\mathbb{R})$ is not complete under this topology.

Let $l, k \in \mathbb{N}$ with $k > l$. The difference between the terms of the sequence is
\[
\phi_{k}(x) - \phi_{l}(x) = \displaystyle\sum_{m=l+1}^{k} \frac{1}{m}\phi(x-m).
\]
Since $\operatorname{supp}(\phi) \subseteq [0,1]$, the support of each translated term $\phi(x-m)$ is contained in the interval $[m, m+1]$. The interiors of these intervals are disjoint and all derivatives of $\phi$ vanish at their endpoints. Therefore, at each point there is at most one nonzero term, even after differentiation.

We evaluate the $j$th seminorm of the difference. For any derivative of order $i \leq j$, we have
\[
|D^{i}(\phi_{k} - \phi_{l})(x)|
=
\left|
\displaystyle\sum_{m=l+1}^{k} \frac{1}{m} D^{i}\phi(x-m)
\right|.
\]
Since at each point there is at most one nonzero summand, the supremum over $\mathbb{R}$ is obtained by maximizing over each interval $[m, m+1]$, giving
\[
\|D^{i}(\phi_{k} - \phi_{l})\|_{L^\infty(\mathbb{R})}
=
\max_{l+1 \leq m \leq k}
\left(
\frac{1}{m}
\|D^{i}\phi\|_{L^\infty(\mathbb{R})}
\right).
\]
Since $\phi\in C_{c}^{\infty}(\mathbb{R})$, its derivatives are bounded, so
\[
C_{j} := \max_{i \leq j} \|D^{i}\phi\|_{L^\infty(\mathbb{R})} < \infty.
\]
Then
\[
\|\phi_{k} - \phi_{l}\|_{j} \leq \frac{1}{l+1} C_{j}.
\]
As $l, k \to \infty$, the quantity $\displaystyle\frac{1}{l+1} C_{j}$ tends to zero, which shows that $(\phi_{k})_{k\in \mathbb{N}}$ is a Cauchy sequence with respect to each seminorm in the family.

In particular, for $j=0$ we obtain
\[
\|\phi_k-\phi_l\|_{L^\infty(\mathbb{R})}
\le
\frac{C_0}{l+1}\to 0,\qquad \text{if }l,k\to\infty,
\]
so $(\phi_k)$ is uniformly Cauchy on $\mathbb{R}$, and consequently there exists a bounded function $\Phi$ such that
\[
\|\phi_k-\Phi\|_{L^\infty(\mathbb{R})}\longrightarrow 0,
\qquad k\to\infty,
\]
that is, $\phi_k\to\Phi$ uniformly on $\mathbb{R}$.

The same argument applied to each derivative shows that for every $i\in\mathbb{N}_0$ the sequence $(D^{i}\phi_k)$ is uniformly Cauchy on $\mathbb{R}$, and therefore there exists a continuous bounded function $\psi_i$ such that
\[
\|D^{i}\phi_k-\psi_i\|_{L^\infty(\mathbb{R})}\longrightarrow 0,
\qquad k\to\infty.
\]
Since also $\phi_k(0)=0$ for every $k$, the sequence converges pointwise on $0$. The theorem on uniform convergence of derivatives implies that $\Phi\in C^\infty(\mathbb{R})$ and that
\[
D^{i}\Phi(x)=\psi_i(x)
\qquad\text{for every }x\in\mathbb{R},\ i\in\mathbb{N}_0.
\]
Moreover, since $\phi_k$ are the partial sums of the series, we have
\[
\Phi(x) = \displaystyle\sum_{m=1}^{\infty} \frac{1}{m}\phi(x-m).
\]

To prove convergence in each seminorm, fix $j\in\mathbb{N}_0$ and $k\in\mathbb{N}$. For every $l>k$, the triangle inequality gives
\[
\|\phi_k-\Phi\|_j
\le
\|\phi_k-\phi_l\|_j
+
\|\phi_l-\Phi\|_j.
\]
By the estimate already obtained, for every $l>k$ we have
\[
\|\phi_k-\phi_l\|_j\le \frac{C_j}{k+1}.
\]
On the other hand, since $D^{i}\phi_l\to D^{i}\Phi$ uniformly on $\mathbb{R}$ for every $i\le j$, we have
\[
\|\phi_l-\Phi\|_j
=
\max_{i\le j}\|D^i\phi_l-D^i\Phi\|_{L^\infty(\mathbb{R})}
\longrightarrow 0,
\qquad l\to\infty.
\]
Letting $l\to\infty$ in the preceding triangle inequality gives
\[
\|\phi_k-\Phi\|_j\le \frac{C_j}{k+1}\to 0,\qquad \text{if }k\to\infty.
\]
Since $j$ was arbitrary, we conclude that for each $j\in\mathbb{N}_0$ we have $\|\phi_k-\Phi\|_j\to 0$ as $k\to\infty$. This shows that $\phi_k$ converges to $\Phi$ with respect to each seminorm in the family, and therefore $\phi_{k}$ converges to $\Phi$ with respect to the topology generated by this family of seminorms.

To analyze the support of $\Phi$, note that by hypothesis $\phi > 0$ on $(0,1)$. For each $m \in \mathbb{N}$, choose the point $x_{m} = m + \displaystyle\frac{1}{2}$. Evaluating the limit function at this point gives
\[
\Phi\left(m + \frac{1}{2}\right)
=
\frac{1}{m}\phi\left(\frac{1}{2}\right) > 0.
\]
This implies that $\Phi$ takes strictly positive values at arbitrarily large points, so its support is unbounded. Consequently, $\operatorname{supp}(\Phi)$ is not compact and $\Phi \notin C_{c}^{\infty}(\mathbb{R})$.
\end{example}
Although the family of seminorms $\{\|\cdot\|_{j}\}_{j\in \mathbb{N}_{0}}$ does not generate the topology of $\mathcal{D}(\Omega)$—since these seminorms do not control the supports of the functions—it is possible to describe the topology $\tau_{LF}$ using suitable families of seminorms that simultaneously incorporate the behavior of the derivatives and of the supports. The following exercise presents one such construction.
\begin{exercise}\label{ejer:distribuciones-y-fourier-abierto-sucesion-compactos-tal-sucesion-creciente}
Let $\Omega\subseteq \mathbb{R}^{n}$ be open. Consider a sequence $\{K_{l}\}_{l\in \mathbb{N}_{0}}$ of compact sets in $\Omega$, with $K_{0}:=\varnothing$, such that for $l\geq 1$ the sequence is increasing and nested, $K_l\subseteq \operatorname{int}(K_{l+1})$, and $\Omega=\displaystyle\bigcup_{l=1}^{\infty}K_{l}$.

Let $\mathbf{m}:=\{m_{l}\}_{l\in \mathbb{N}_{0}}$ be a sequence in $\mathbb{N}_{0}$ such that $m_{l}\to \infty$, and let $\mathbf{a}:=\{a_{l}\}_{l\in \mathbb{N}_{0}}$ be a sequence in $(0,\infty)$ such that $a_{l}\to 0$. For each $\phi\in \mathcal{D}(\Omega)$ we define
\[
p_{\mathbf{m},\mathbf{a}}(\phi)
:=
\sup_{l\in \mathbb{N}_{0}}
\sup_{x\in \Omega\setminus K_{l}}
\frac{1}{a_{l}}
\displaystyle\sum_{|\alpha|\leq m_{l}}
|D^{\alpha}\phi(x)|.
\]

\begin{enumerate}
\item Prove that $p_{\mathbf{m},\mathbf{a}}$ is a seminorm on $\mathcal{D}(\Omega)$.
\item Prove that the family of seminorms
$\{p_{\mathbf{m},\mathbf{a}}\}_{\mathbf{m},\mathbf{a}}$
generates the topology $\tau_{LF}$ of Theorem~\ref{teo: topologia de D(Omega)}.
\end{enumerate}
\end{exercise}
Observe that the preceding family of seminorms is, in general, uncountable. This reflects the fact that $\mathcal{D}(\Omega)$ will turn out not to be metrizable (as the reader will justify later) and, in particular, is not a Fréchet space.

Fixing the support recovers exactly the topology we already defined on $\mathcal D_K(\Omega)$.
\begin{theorem}\label{teo: Dk tiene topologia de subespacio en D}\index{test function space!subspace topology of \(\mathcal D_K(\Omega)\)}
 Let $\Omega\subseteq \mathbb{R}^{n}$ be open. For any compact set $K\subseteq \Omega$, the topology $\tau_{K}$ agrees with the subspace topology on $\mathcal{D}_{K}(\Omega)$ induced by $\tau_{LF}$ as a subset of $\mathcal{D}(\Omega)$.
\end{theorem}
\begin{proof}
The inclusion $\mathcal D_K(\Omega)\hookrightarrow\mathcal D(\Omega)$ is continuous: the inverse image of each $U\in\mathcal B_0$ is, by definition, a neighborhood of the origin in $\mathcal D_K(\Omega)$. Thus the induced topology is no finer than $\tau_K$.

Conversely, for $j\in\mathbb N_0$ and $\varepsilon>0$, the ball
\[
V_{j,\varepsilon}:=\{\phi\in C_c^\infty(\Omega)
\mid \|\phi\|_j<\varepsilon\}
\]
belongs to $\mathcal B_0$: on any compact support, the global norm in \eqref{eq: normas globales en D(Omega)} agrees with the norm corresponding to that compact set. Hence
\[
V_{j,\varepsilon}\cap\mathcal D_K(\Omega)
=\{\phi\in\mathcal D_K(\Omega)\mid\|\phi\|_{C^j(K)}<\varepsilon\}
\]
is a neighborhood for the induced topology. These balls form a neighborhood base for $\tau_K$, proving the reverse inclusion.
\end{proof}
We next study the topologically bounded sets\footnote{Let $X$ be a topological vector space. A subset $B\subseteq X$ is said to be topologically bounded if for every neighborhood $U$ of the origin there exists $\lambda>0$ such that $B\subseteq \lambda U$.} in $\mathcal{D}(\Omega)$, for which we will need the following result:
\begin{lemma}\label{lem: evaluacion abierta en D_K}
Let $\Omega\subseteq \mathbb{R}^{n}$ be an open set, let $K\subseteq \Omega$ be compact, and let $x_{0}\in \Omega$. Then, for every $r>0$, the set
$U:=\{\phi\in \mathcal{D}_{K}(\Omega)\mid |\phi(x_{0})|<r\}$
is open in $(\mathcal{D}_{K}(\Omega),\tau_{K})$.
\end{lemma}

\begin{proof}
First observe that if $x_{0}\notin K$, then every function $\phi\in \mathcal{D}_{K}(\Omega)$ satisfies $\phi(x_{0})=0$, so $U=\mathcal{D}_{K}(\Omega)$, which is open.

Now suppose that $x_{0}\in K$. Let $\phi\in U$, that is, $|\phi(x_{0})|<r$. Define $\varepsilon:=r-|\phi(x_{0})|>0$. Consider the neighborhood of $\phi$ given by $V:=\{\psi\in \mathcal{D}_{K}(\Omega)\mid \|\psi-\phi\|_{C^0(K)}<\varepsilon\}$. Then, for every $\psi\in V$, we have in particular $|\psi(x_{0})-\phi(x_{0})|\leq \|\psi-\phi\|_{C^{0}(K)}<\varepsilon$.

Therefore, $\displaystyle |\psi(x_{0})|\leq |\phi(x_{0})|+|\psi(x_{0})-\phi(x_{0})|<|\phi(x_{0})|+\varepsilon=r$. This shows that $\psi\in U$, and consequently $V\subseteq U$.

We have proved that every $\phi\in U$ has a neighborhood contained in $U$, so $U$ is open in $\tau_{K}$.
\end{proof}
\begin{theorem}\label{teo: conjuntos topologicamente acotados en D(Omega)}\index{topologically bounded sets in D(Omega)}
 Let $\Omega\subseteq \mathbb{R}^{n}$ be an open set and let $W\subseteq \mathcal{D}(\Omega)$ be a topologically bounded set. Then there exists a compact set $K\subseteq \Omega$ such that $W\subseteq \mathcal{D}_{K}(\Omega)$. Moreover, for each $j\in \mathbb{N}_{0}$, there exists $M_{j}>0$ such that $\|\phi\|_{j}\leq M_{j}$ for every $\phi\in W$, where $\|\cdot\|_{j}$ is given by \eqref{eq: normas globales en D(Omega)}.
\end{theorem}

\begin{proof}
We argue by contradiction. Suppose that $W\not\subseteq \mathcal{D}_{K}(\Omega)$ for every compact set $K\subseteq \Omega$. Let $\{K_l\}_{l\in\mathbb{N}}$ be an increasing sequence of compact sets such that $K_l\subseteq \operatorname{int}(K_{l+1})$ and $\Omega=\displaystyle\bigcup_{l=1}^{\infty}K_l$.

Then, for each $l\in\mathbb{N}$, there exists $\phi_l\in W$ such that $\supp(\phi_l)\not\subseteq K_l$. Therefore, we can choose a point $y_l\in \supp(\phi_l)\setminus K_l$.

Since $K_l$ is closed in $\Omega$, the set $\Omega\setminus K_l$ is an open neighborhood of $y_l$. Moreover, by definition of support, $y_l\in\supp(\phi_l)$ implies that every neighborhood of $y_l$ intersects the set $\{x\in\Omega\mid \phi_l(x)\neq 0\}$. In particular, there exists a point $x_l\in \Omega\setminus K_l$ such that $\phi_l(x_l)\neq 0$.

We define
\[
U:=\left\{\phi\in \mathcal{D}(\Omega)\Biggm|
|\phi(x_l)|<\frac{|\phi_l(x_l)|}{l},\ \forall l\in\mathbb{N}
\right\}.
\]

We claim that $U\in\mathcal{B}_0$. First, $0\in U$, since $|0(x_l)|=0<\displaystyle\frac{|\phi_l(x_l)|}{l}$ for every $l\in\mathbb{N}$.

Moreover, $U$ is balanced. If $\phi\in U$ and $|\lambda|\leq 1$, then
\[
|(\lambda\phi)(x_l)|
=
|\lambda|\,|\phi(x_l)|
\leq |\phi(x_l)|
<
\frac{|\phi_l(x_l)|}{l}
\]
for every $l\in\mathbb{N}$. Therefore, $\lambda\phi\in U$.

Also, $U$ is convex. If $\phi,\psi\in U$ and $s\in[0,1]$, then, for every $l\in\mathbb{N}$,
\[
|(s\phi+(1-s)\psi)(x_l)|
\leq
s|\phi(x_l)|+(1-s)|\psi(x_l)|
<
s\frac{|\phi_l(x_l)|}{l}+(1-s)\frac{|\phi_l(x_l)|}{l}
=
\frac{|\phi_l(x_l)|}{l}.
\]
Thus, $s\phi+(1-s)\psi\in U$.

It remains to check that, for every compact set $K\subseteq\Omega$, the set $U\cap\mathcal{D}_K(\Omega)$ is a neighborhood of the origin in $(\mathcal{D}_K(\Omega),\tau_K)$. Let $K\subseteq\Omega$ be compact. Since $\Omega=\displaystyle\bigcup_{l=1}^{\infty}K_l$ and the sequence $(K_l)$ is increasing, there exists $L\in\mathbb{N}$ such that $K\subseteq K_L$. Then, if $l\geq L$, we have $K\subseteq K_l$ and, since $x_l\notin K_l$, we conclude that $x_l\notin K$.

We define
\[
I_K:=\{l\in\mathbb{N}\mid x_l\in K\}.
\]
By the preceding argument, $I_K\subseteq\{1,\dots,L-1\}$, so $I_K$ is finite.

If $\phi\in\mathcal{D}_K(\Omega)$ and $l\notin I_K$, then $x_l\notin K$ and, since $\supp(\phi)\subseteq K$, we have $\phi(x_l)=0$. Consequently, the conditions corresponding to these indices are automatically satisfied. Therefore,
\[
U\cap \mathcal{D}_{K}(\Omega)
=
\left\{\phi\in \mathcal{D}_{K}(\Omega)\Biggm|
|\phi(x_l)|<\frac{|\phi_l(x_l)|}{l},\ \forall l\in I_K
\right\}.
\]

Since $I_K$ is finite, this set is a finite intersection of sets of the form
\[
\left\{\phi\in\mathcal{D}_K(\Omega)\mid |\phi(x_l)|<c_l\right\},
\]
which are open in $(\mathcal{D}_K(\Omega),\tau_K)$ by Lemma~\ref{lem: evaluacion abierta en D_K}. Moreover, they contain the origin. Consequently, $U\cap\mathcal{D}_K(\Omega)$ is a neighborhood of the origin in $(\mathcal{D}_K(\Omega),\tau_K)$.

Since this holds for every compact set $K\subseteq\Omega$, we conclude that $U\in\mathcal{B}_0$. In particular, $U=0+U\in\mathcal{B}$, so $U$ is a neighborhood of the origin in $(\mathcal{D}(\Omega),\tau_{LF})$.

Since $W$ is topologically bounded, there exists $t>0$ such that $W\subseteq tU$. Let $l\in\mathbb{N}$ satisfy $l\geq t$. Then
\[
\frac{|\phi_l(x_l)|}{t}\geq \frac{|\phi_l(x_l)|}{l},
\]
and therefore
\[
\left|t^{-1}\phi_l(x_l)\right|
=
\frac{|\phi_l(x_l)|}{t}
\geq
\frac{|\phi_l(x_l)|}{l}.
\]
Thus, $t^{-1}\phi_l\notin U$, which is equivalent to $\phi_l\notin tU$. This contradicts the fact that $\phi_l\in W\subseteq tU$.

Therefore, there exists a compact set $K\subseteq\Omega$ such that $W\subseteq\mathcal{D}_K(\Omega)$.

Finally, by Theorem~\ref{teo: Dk tiene topologia de subespacio en D}, the set $W=W\cap \mathcal{D}_K(\Omega)$ is bounded with respect to the topology $\tau_K$. By Corollary~\ref{cor: topologia generada por seminormas hausdorff y topologicamente acotado}, for each $j\in\mathbb{N}_0$ there exists $M_j>0$ such that
\[
\|\phi\|_{C^{j}(K)}\leq M_j
\qquad
\text{for every } \phi\in W.
\]
Since $\supp(\phi)\subseteq K$ for every $\phi\in W$, we have $\|\phi\|_j=\|\phi\|_{C^{j}(K)}$ on $W$. Therefore,
\[
\|\phi\|_j\leq M_j
\qquad
\text{for every } \phi\in W.
\]
This completes the proof.
\end{proof}
The space $\mathcal{D}(\Omega)$ with the topology $\tau_{LF}$ of Theorem~\ref{teo: topologia de D(Omega)} will turn out to be a (sequentially) complete topological vector space, and we will be able to characterize convergence of sequences quite conveniently:
\begin{theorem}\label{teo: convergencia en D(Omega)}\index{convergence in D(Omega)}
Let $\Omega\subseteq \mathbb{R}^{n}$ be an open set. Then the space $(\mathcal{D}(\Omega),\tau_{LF})$ is sequentially complete. Moreover, a sequence $\{\phi_{k}\}_{k\in \mathbb{N}}\subseteq \mathcal{D}(\Omega)$ converges to $\phi\in \mathcal{D}(\Omega)$ with respect to $\tau_{LF}$ if and only if:
\begin{enumerate}
\item there exists a compact set $K\subseteq \Omega$ such that $\supp(\phi_{k})\subseteq K$ for every $k\in \mathbb{N}$,
\item $\displaystyle\lim_{k\to \infty}D^{\alpha}\phi_{k}=D^{\alpha}\phi$ uniformly on $\Omega$ for every multi-index $\alpha\in \mathbb{N}_{0}^{n}$.
\end{enumerate}
\end{theorem}

\begin{proof}
Let $\{\phi_k\}_{k\in\mathbb N}$ be a Cauchy sequence in $\mathcal D(\Omega)$. By Proposition~\ref{prop: sucesion de Cauchy tiene rango topologicamente acotado}, the set $\{\phi_k\mid k\in\mathbb N\}$ is topologically bounded in $\mathcal D(\Omega)$. By Theorem~\ref{teo: conjuntos topologicamente acotados en D(Omega)}, there exists a compact set $K\subseteq\Omega$ such that
\[
\{\phi_k\mid k\in\mathbb N\}\subseteq\mathcal D_K(\Omega).
\]

By Theorem~\ref{teo: Dk tiene topologia de subespacio en D}, the sequence $\{\phi_k\}_{k\in\mathbb N}$ is Cauchy in $\mathcal D_K(\Omega)$ with respect to $\tau_K$. Since $\mathcal D_K(\Omega)$ is complete, there exists $\phi\in\mathcal D_K(\Omega)$ such that
\[
\phi_k\longrightarrow\phi
\]
in $\mathcal D_K(\Omega)$ with respect to $\tau_K$.

Applying Theorem~\ref{teo: Dk tiene topologia de subespacio en D} again, we conclude that $\phi_k\to\phi$ in $\mathcal D(\Omega)$ with respect to $\tau_{LF}$. Therefore, $\mathcal D(\Omega)$ is sequentially complete.

We now prove the equivalence. Let $\{\phi_k\}_{k\in\mathbb{N}}\subseteq \mathcal{D}(\Omega)$ and let $\phi\in\mathcal{D}(\Omega)$.

$(\Rightarrow)$ Suppose that $\phi_k\to \phi$ with respect to $\tau_{LF}$. Then the sequence is Cauchy and, by Proposition~\ref{prop: sucesion de Cauchy tiene rango topologicamente acotado}, $\{\phi_k\mid k\in\mathbb{N}\}$ is topologically bounded in $\mathcal{D}(\Omega)$. By Theorem~\ref{teo: conjuntos topologicamente acotados en D(Omega)}, there exists a compact set $K\subseteq \Omega$ such that
\[
\phi_k\in \mathcal{D}_K(\Omega)
\]
for every $k\in \mathbb{N}$. That is, $\supp(\phi_k)\subseteq K$ for every $k\in \mathbb{N}$. Since also $\phi\in\mathcal{D}(\Omega)$, we have $\phi\in \mathcal{D}_{K'}(\Omega)$ for some compact set $K'\subseteq\Omega$. Replacing $K$ by $K\cup K'$, we may assume that $\phi\in\mathcal{D}_K(\Omega)$.

By Theorem~\ref{teo: Dk tiene topologia de subespacio en D}, the topology induced by $\tau_{LF}$ on $\mathcal{D}_K(\Omega)$ agrees with $\tau_K$. In particular, $\phi_k\to\phi$ in $\mathcal{D}_K(\Omega)$ with respect to $\tau_K$. By definition of $\tau_K$, this is equivalent to saying that, for every $j\in\mathbb{N}_0$,
\[
\|\phi_k-\phi\|_{C^{j}(K)}\to 0.
\]
Consequently, for every multi-index $\alpha\in\mathbb{N}_{0}^{n}$ we have
\[
D^\alpha\phi_k\to D^\alpha\phi
\]
uniformly on $K$. Since $\phi_k,\phi\in\mathcal D_K(\Omega)$, their derivatives vanish on $\Omega\setminus K$. Consequently,
\[
\|D^\alpha\phi_k-D^\alpha\phi\|_{L^\infty(\Omega)}
=
\|D^\alpha\phi_k-D^\alpha\phi\|_{C^0(K)}
\longrightarrow0.
\]
The convergence is therefore uniform on $\Omega$. This proves both conditions.

$(\Leftarrow)$ Now suppose that there exists a compact set $K\subseteq\Omega$ such that $\supp(\phi_k)\subseteq K$ for every $k\in\mathbb{N}$ and that
\[
D^\alpha\phi_k\to D^\alpha\phi
\]
uniformly on $\Omega$ for every multi-index $\alpha\in\mathbb{N}_{0}^{n}$. The case $\alpha=0$ implies pointwise convergence throughout $\Omega$. If $x\in\Omega\setminus K$, then $\phi_k(x)=0$ for every $k$, so $\phi(x)=\lim_{k\to\infty}\phi_k(x)=0$. Since $K$ is closed in $\Omega$, it follows that $\supp(\phi)\subseteq K$. Then $\phi_k,\phi\in\mathcal{D}_K(\Omega)$ and, for every $j\in\mathbb{N}_0$,
\[
\|\phi_k-\phi\|_{C^{j}(K)}
=
\max_{|\alpha|\le j}\|D^\alpha\phi_k-D^\alpha\phi\|_{C^0(K)}
\to 0.
\]
Therefore, $\phi_k\to\phi$ in $\mathcal{D}_K(\Omega)$ with respect to $\tau_K$. Since $\tau_K$ agrees with the subspace topology induced by $\tau_{LF}$ on $\mathcal{D}_K(\Omega)$, we conclude that $\phi_k\to\phi$ in $\mathcal{D}(\Omega)$ with respect to $\tau_{LF}$.
\end{proof}
When $\Omega\subseteq\mathbb R^n$ is open and nonempty, with $n\geq1$, the space $\mathcal{D}(\Omega)$ is not metrizable:
\begin{exercise}\label{ejer:distribuciones-y-fourier-abierto-pruebe-compacto-cerrado-interior-vacio}
 Let $\Omega\subseteq \mathbb{R}^{n}$ be open and nonempty, with $n\geq1$.
 \begin{enumerate}
 \item Prove that for every compact set $K\subseteq \Omega$, $\mathcal{D}_{K}(\Omega)$ is closed and has empty interior in $\mathcal{D}(\Omega)$.
 \item Prove that $\mathcal{D}(\Omega)$ is not metrizable.
 \end{enumerate}
 \textit{Hint:} To prove that $\mathcal{D}_{K}(\Omega)$ has empty interior in the first part, construct a sequence of functions in $\mathcal{D}(\Omega) \setminus \mathcal{D}_{K}(\Omega)$ converging to the constant $0$ in the topology of $\mathcal{D}(\Omega)$. For the second part, consider an exhaustion by compact sets $(K_{j})_{j\in \mathbb{N}}$ such that $\Omega = \displaystyle\bigcup_{j=1}^{\infty} K_{j}$. Use the preceding part and the Baire category theorem (Theorem~\ref{teo:categoria-de-baire}), recalling the relation between complete metric spaces and Baire spaces.
\end{exercise}
\section{LF spaces}
The spaces $\mathcal{D}(\Omega)$ can be obtained as a particular case of a general construction in functional analysis known as a \textit{strict inductive limit of Fréchet spaces}. Spaces obtained in this way are called \textit{LF spaces}.

To formalize this notion, recall that, given a family of maps $g_i\colon X_i\longrightarrow X$, there is a final topology on $X$, characterized as the finest topology making all the $g_i$ continuous. In the setting of topological vector spaces, however, this purely topological final topology need not be compatible with the linear structure of $X$, nor need it be locally convex.

Thus, when the spaces $X_i$ are locally convex and the maps $g_i$ are linear, we are interested in the finest among the locally convex topologies on $X$ making all the $g_i$ continuous. In the case of inclusions of subspaces, this construction takes the following form.
\begin{theorem}[Locally convex final topology]\label{teo: topologia-final-lc}\index{locally convex final topology}
Let $X$ be a vector space and let $(X_i)_{i\in I}$ be a family of subspaces of $X$ such that
\[
X=\bigcup_{i\in I}X_i.
\]
Suppose that each $X_i$ is equipped with a locally convex topology $\tau_i$. Then there is a unique locally convex topology $\tau_{LC}$ on $X$ that is the finest among all locally convex topologies making the canonical inclusions
\[
\iota_i\colon X_i\hookrightarrow X,
\qquad i\in I.
\]
continuous.

The topology $\tau_{LC}$ is called the \textit{locally convex final topology} associated with the family $\bigl((X_i,\tau_i),\iota_i\bigr)_{i\in I}$.
\end{theorem}
\begin{proof}
Let $\mathcal{U}$ be the collection of all absolutely convex subsets\footnotemark $U\subseteq X$ such that, for each $i\in I$, $U\cap X_i$ is a neighborhood of the origin in $(X_i,\tau_i)$.

First observe that every element of $\mathcal{U}$ is absorbing. Let $U\in\mathcal{U}$ and let $x\in X$. Since
\[
X=\bigcup_{i\in I}X_i,
\]
there exists $i\in I$ such that $x\in X_i$. By hypothesis, $U\cap X_i$ is a neighborhood of the origin in the topological vector space $X_i$ and is therefore absorbing in $X_i$. Consequently, there exists $\lambda>0$ such that
\[
x\in t(U\cap X_i)
\qquad
\text{for every } t\geq \lambda.
\]
Since $t(U\cap X_i)\subseteq tU$ for every $t>0$, it follows that
\[
x\in tU
\qquad
\text{for every } t\geq \lambda.
\]
Since $x\in X$ was arbitrary, we conclude that $U$ is absorbing.

We now show that $\mathcal{U}$ satisfies the hypotheses of Proposition~\ref{prop: base-discos-topologia-lc}.

\begin{enumerate}
 \item By definition, each $U\in\mathcal{U}$ is absolutely convex, and we have just proved that it is also absorbing.

 \item Let $U,V\in\mathcal{U}$. Then $U\cap V$ is absolutely convex. Moreover, for each $i\in I$,
 \[
 (U\cap V)\cap X_i
 =
 (U\cap X_i)\cap (V\cap X_i).
 \]
 Since $U\cap X_i$ and $V\cap X_i$ are neighborhoods of the origin in $(X_i,\tau_i)$, their intersection is as well. Therefore, $U\cap V\in\mathcal{U}$.

 \item Let $U\in\mathcal{U}$ and let $\lambda>0$. Define $W:=\displaystyle\frac{\lambda}{2}U$. Since $U$ is absolutely convex, so is $W$. Moreover, for each $i\in I$, since $X_i$ is a subspace of $X$, we have
\[
W\cap X_i
=
\left(\frac{\lambda}{2}U\right)\cap X_i
=
\frac{\lambda}{2}(U\cap X_i).
\]
Since multiplication by the nonzero scalar $\displaystyle\frac{\lambda}{2}$ is a homeomorphism of $(X_i,\tau_i)$ onto itself, the set $\displaystyle\frac{\lambda}{2}(U\cap X_i)$ is a neighborhood of the origin in $X_i$. Thus, $W$ is absolutely convex and $W\cap X_i$ is a neighborhood of the origin in $X_i$ for every $i\in I$. By definition of $\mathcal{U}$, we conclude that $W\in\mathcal{U}$.

Finally, since $U$ is balanced and $\displaystyle\frac{\lambda}{2}\leq \lambda$, we have
\[
W=\frac{\lambda}{2}U\subseteq \lambda U.
\]
\end{enumerate}

By Proposition~\ref{prop: base-discos-topologia-lc}, the family $\mathcal{U}$ is a local base at the origin for a unique locally convex topology $\tau_{LC}$ on $X$.

We show that each canonical inclusion
\[
\iota_i\colon X_i\hookrightarrow X
\]
is continuous with respect to $\tau_{LC}$. Since $\iota_i$ is linear, it suffices to check continuity at the origin. Let $A$ be a neighborhood of the origin in $(X,\tau_{LC})$. Since $\mathcal{U}$ is a local base at the origin, there exists $U\in\mathcal{U}$ such that $U\subseteq A$. Then $\iota_i^{-1}(U)=U\cap X_i$ is a neighborhood of the origin in $(X_i,\tau_i)$, and $\iota_i^{-1}(U)\subseteq \iota_i^{-1}(A)$. Therefore, $\iota_i^{-1}(A)$ is a neighborhood of the origin in $X_i$, proving that $\iota_i$ is continuous.

Finally, we show that $\tau_{LC}$ is the finest among all locally convex topologies with this property. Let $\tau'$ be another locally convex topology on $X$ such that all the inclusions $\iota_i$ are continuous. Let $V$ be an absolutely convex neighborhood of the origin in $(X,\tau')$. Then, for each $i\in I$,
\[
V\cap X_i=\iota_i^{-1}(V)
\]
is a neighborhood of the origin in $(X_i,\tau_i)$, by continuity of $\iota_i$. Therefore, $V\in\mathcal{U}$.

Now, by Proposition~\ref{prop: localmente convexo tiene una base de vecindades absolutamente convexas}, the absolutely convex neighborhoods of the origin form a local base in $(X,\tau')$. Thus, if $A$ is a neighborhood of the origin in $(X,\tau')$, there is an absolutely convex neighborhood $V$ of the origin such that $V\subseteq A$. By the preceding argument, $V\in\mathcal{U}$. Since $\mathcal{U}$ is a local base at the origin for $\tau_{LC}$, the set $V$ is a neighborhood of the origin with respect to $\tau_{LC}$; since it is contained in $A$, we conclude that $A$ is also a neighborhood of the origin with respect to $\tau_{LC}$. We have proved that every neighborhood of the origin for $\tau'$ is also a neighborhood of the origin for $\tau_{LC}$. Consequently, $\tau'\subseteq\tau_{LC}$.

Thus, $\tau_{LC}$ is the finest among all locally convex topologies making the inclusions $\iota_i$ continuous.

Uniqueness follows immediately from this maximality property.
\end{proof}

\footnotetext{A subset $U$ of a vector space is said to be \textit{absolutely convex} (or a \textit{disk}) if it is both convex and balanced.}

With this construction, we can formally define LF spaces.

\begin{definition}\label{def: espacio LF}\index{strict inductive limit}
 Let $X$ be a vector space over $\mathbb{R}$ or $\mathbb{C}$. Suppose there is a sequence of subspaces $(X_{n})_{n\in \mathbb{N}}$ such that:
 \begin{enumerate}
 \item $X_{1} \subseteq X_{2} \subseteq X_{3} \subseteq \dots$ and $X=\displaystyle\bigcup_{n=1}^{\infty}X_{n}$.
 \item For each $n$, the space $X_{n}$, equipped with a topology $\tau_{n}$, is a Fréchet space.
 \item For each $n$, the topology $\tau_{n}$ agrees with the subspace topology induced on $X_{n}$ by $(X_{n+1},\tau_{n+1})$.
 \end{enumerate}
 Let $\tau_{LF}$ be the locally convex final topology on $X$ associated with the inclusions $\iota_n\colon (X_n,\tau_n)\hookrightarrow X$, guaranteed by Theorem~\ref{teo: topologia-final-lc}. Then $(X,\tau_{LF})$ is called the \textit{strict inductive limit} of the sequence $(X_n,\tau_n)$, and we write
 \[ (X,\tau_{LF})=\varinjlim_{n\in\mathbb{N}} (X_n,\tau_n).
 \]
 In this case, $(X,\tau_{LF})$ is said to be an \textit{LF space}. The adjective \textit{strict} refers to compatibility of the topologies of the stages, not to the inclusions being proper. Stationary sequences are allowed; in particular, every Fréchet space is obtained by taking $X_n=X$ with the same topology for every $n$.
\end{definition}
\begin{remark}\label{obs:distribuciones-y-fourier-definicion-anterior-espacio-banach-denomina-espacio}
If in the preceding definition each space $X_n$ is Banach, then $X$ is called an LB space. In particular, every LB space is an LF space.
\end{remark}
\begin{remark}\label{obs:distribuciones-y-fourier-nocion-espacio-lf-surge-naturalmente-desarrollo}
The notion of an LF space arose naturally in the development of functional analysis and distribution theory during the first half of the twentieth century. At that time, it was observed that many fundamental spaces in analysis could not be adequately described within the category of Banach spaces or even Fréchet spaces, although they could be constructed as increasing unions of the latter.

One of the most important examples appeared in Laurent Schwartz's distribution theory. The space of smooth functions with compact support,
$
\mathcal D(\Omega)=C_c^\infty(\Omega),
$
has a topological structure essential to the study of distributions and partial differential equations, but is not metrizable under its natural topology. However, it can be written as a countable increasing union of the spaces
$
\mathcal D_K(\Omega)
=
\{f\in C^\infty(\Omega)\mid \operatorname{supp}(f)\subseteq K\},
$
where $K\Subset\Omega$ ranges over the compact subsets of $\Omega$. Each of these spaces turns out to be a Fréchet space, which motivated the systematic study of countable inductive limits of Fréchet spaces.

The theory was subsequently developed within the general framework of locally convex spaces by authors such as Jean Dieudonné and Alexander Grothendieck, becoming a fundamental tool in functional analysis, distribution theory, differential geometry, and global analysis.
\end{remark}
To show that $\mathcal{D}(\Omega)$ is an LF space, we first construct an exhaustion of $\Omega$ by compact sets.

For each $l\in\mathbb N$, we define
$
C_l
:=
\left\{
x\in\Omega
\ \middle|\
d\bigl(x,\mathbb R^n\setminus\Omega\bigr)\geq \displaystyle\frac{1}{l}
\right\},
$
and then
$
K_l:=C_l\cap \overline{B}_{\mathrm{euc}}(0,l).
$

The geometric idea behind this construction is twofold: on the one hand, the sets $C_l$ force points to remain at a positive distance from the boundary $\partial\Omega$; on the other, intersection with the closed ball $\overline{B}_{\mathrm{euc}}(0,l)$ prevents points from escaping to infinity.

It is immediate to check that each $K_l$ is compact, that
$
K_l\subseteq \operatorname{int}(K_{l+1}),
$
and that
$
\Omega=\displaystyle\bigcup_{l=1}^{\infty}K_l.
$

\begin{figura}[H]

    \begin{minipage}[c]{0.53\textwidth}
        \centering
        \includegraphics[width=\textwidth]{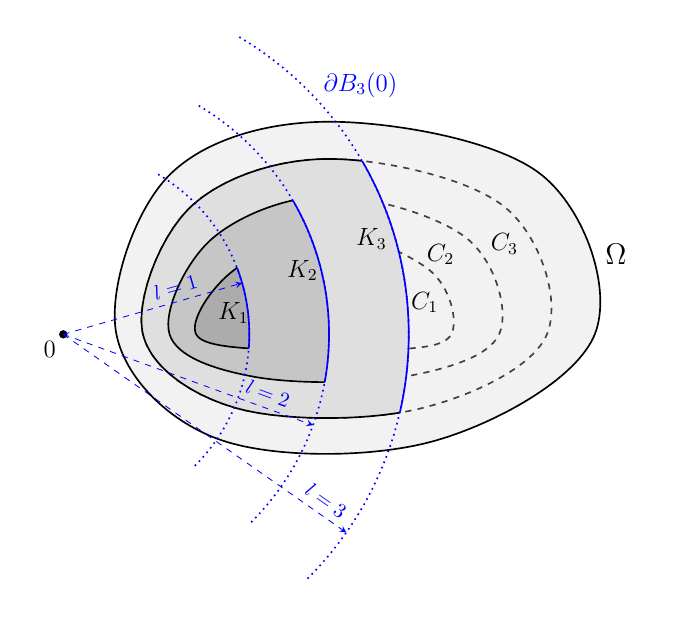}
    \end{minipage}
    \hfill
    \begin{minipage}[c]{0.41\textwidth}
        \small

        \begin{center}
            \textbf{Geometric interpretation}
        \end{center}

        \vspace{0.3em}

        The set
        $
        K_l=C_l\cap \overline{B}_{\mathrm{euc}}(0,l)
        $
        is obtained by combining two complementary restrictions:

        \begin{enumerate}[label=\textup{(\alph*)}, leftmargin=*]
            \item $C_l$ removes points too close to $\partial\Omega$;

            \item $\overline{B}_{\mathrm{euc}}(0,l)$ controls behavior at infinity.
        \end{enumerate}

        Consequently, the boundary of $K_l$ consists of:
        \[
        \partial K_l
        \subseteq
        \Bigl\{
        x\in\Omega
        \ \mid\
        d(x,\partial\Omega)=\displaystyle\frac{1}{l}
        \Bigr\}
        \cup
        \partial B_{\mathrm{euc}}(0,l).
        \]

        In this way, the sequence $(K_l)_{l\in\mathbb N}$ provides a natural exhaustion of $\Omega$ by compact sets.
    \end{minipage}

    \caption{Construction of a compact exhaustion of $\Omega$.}
    \label{fig:exhaustion-compactos-bolas}

\end{figura}

Now taking
$
X_l:=\mathcal D_{K_l}(\Omega),
$
conditions \textup{(a)}--\textup{(c)} of Definition~\ref{def: espacio LF} are satisfied. The family is increasing and exhaustive, each space $X_l$ is Fréchet, and the topologies are compatible by Theorem~\ref{teo: Dk tiene topologia de subespacio en D}.

It remains to check that the topology $\tau_{LF}$ agrees with the locally convex inductive topology associated with the inclusions $\mathcal D_{K_l}(\Omega)\hookrightarrow \mathcal D(\Omega)$.

\begin{proposition}\label{prop: D es limite inductivo}
The topology $\tau_{LF}$ on $\mathcal{D}(\Omega)$ agrees with the locally convex inductive topology associated with the family of inclusions
\[
\iota_l\colon \mathcal{D}_{K_l}(\Omega)\hookrightarrow \mathcal{D}(\Omega),
\qquad l\in\mathbb N.
\]
\end{proposition}

\begin{proof}
Denote by $\tau_{LC}$ the locally convex inductive topology on $\mathcal{D}(\Omega)$ associated with the family
\[
\bigl((\mathcal{D}_{K_l}(\Omega),\tau_{K_l}),\iota_l\bigr)_{l\in\mathbb N}.
\]
By Theorem~\ref{teo: topologia-final-lc}, a local base at the origin for $\tau_{LC}$ is given by the family $\mathcal{U}$ of all absolutely convex subsets $U\subseteq\mathcal{D}(\Omega)$ such that $U\cap \mathcal{D}_{K_l}(\Omega)$ is a neighborhood of the origin in $(\mathcal{D}_{K_l}(\Omega),\tau_{K_l})$ for every $l\in\mathbb N$.

We compare this family with the local base $\mathcal{B}_0$ defining the topology $\tau_{LF}$. Recall that, by definition, $V\in\mathcal{B}_0$ if and only if $0\in V$, $V$ is absolutely convex, and $V\cap \mathcal{D}_{K}(\Omega)$ is a neighborhood of the origin in $(\mathcal{D}_{K}(\Omega),\tau_{K})$ for every compact set $K\subseteq\Omega$.

Thus, the only apparent difference between the families $\mathcal{U}$ and $\mathcal{B}_0$ is that, in the definition of $\mathcal{U}$, the neighborhood condition is required only for the compact sets in the exhaustion $(K_l)_{l\in\mathbb N}$, whereas in the definition of $\mathcal{B}_0$ it is required for every compact set $K\subseteq\Omega$. We show that the two conditions are equivalent.

First we prove that $\mathcal{B}_0\subseteq\mathcal{U}$. Let $V\in\mathcal{B}_0$. Since each $K_l$ is compact, the definition of $\mathcal{B}_0$ implies that $V\cap \mathcal{D}_{K_l}(\Omega)$ is a neighborhood of the origin in $(\mathcal{D}_{K_l}(\Omega),\tau_{K_l})$ for every $l\in\mathbb N$. Since $V$ is also absolutely convex, we conclude that $V\in\mathcal{U}$.

We now prove the reverse inclusion. Let $U\in\mathcal{U}$ and let $K\subseteq\Omega$ be an arbitrary compact set. Since $(K_l)_{l\in\mathbb N}$ is an exhaustion of $\Omega$ by compact sets, there exists $l\in\mathbb N$ such that $K\subseteq K_l$.

Then $\mathcal{D}_{K}(\Omega)\subseteq \mathcal{D}_{K_l}(\Omega)$.

Moreover, the topology $\tau_K$ agrees with the topology induced by $\tau_{K_l}$ on the subspace $\mathcal{D}_{K}(\Omega)$. Since $U\cap \mathcal{D}_{K_l}(\Omega)$ is a neighborhood of the origin in $(\mathcal{D}_{K_l}(\Omega),\tau_{K_l})$, it follows that
\[
\bigl(U\cap \mathcal{D}_{K_l}(\Omega)\bigr)\cap \mathcal{D}_{K}(\Omega)
=
U\cap \mathcal{D}_{K}(\Omega)
\]
is a neighborhood of the origin in $(\mathcal{D}_{K}(\Omega),\tau_K)$.

Since $K\subseteq\Omega$ was an arbitrary compact set, $U$ satisfies the neighborhood condition required in the definition of $\mathcal{B}_0$. Moreover, since $U\cap \mathcal{D}_{K_l}(\Omega)$ is a neighborhood of the origin for any $l\in\mathbb N$, we have $0\in U$. Therefore, $U\in\mathcal{B}_0$.

We have proved that $\mathcal{B}_0=\mathcal{U}$.

Consequently, $\tau_{LF}$ and $\tau_{LC}$ have the same local base at the origin. By Proposition~\ref{prop: topologia-determinada-por-base-local}, the two topologies agree. That is, $\tau_{LF}=\tau_{LC}$.
\end{proof}

Consequently, $(\mathcal{D}(\Omega),\tau_{LF})$ is the strict inductive limit
\[
\mathcal{D}(\Omega)=\varinjlim_{l\in\mathbb{N}}\mathcal{D}_{K_l}(\Omega),
\]
and is therefore an LF space.

The advantage of this description is that many fundamental properties of $\mathcal{D}(\Omega)$ follow directly from its LF structure in a purely topological framework. These properties will allow us to construct spaces of distributions on vector bundles without repeating the constructions and proofs.
\begin{theorem}\label{teo: propiedades LF}\index{LF properties}
Let $(X,\tau)$ be an $LF$ space with defining sequence $(X_n)_{n\in\mathbb N}$. Then:
\begin{enumerate}
 \item $X$ is Hausdorff.

 \item A set $B\subseteq X$ is topologically bounded if and only if there exists $m\in\mathbb N$ such that $B\subseteq X_m$ and $B$ is topologically bounded in $X_m$.

 \item $X$ is complete; in particular, it is sequentially complete.

 \item A sequence $(x_k)_{k\in\mathbb N}$ converges to $x\in X$ if and only if there exists $m\in\mathbb N$ such that $\{x\}\cup\{x_k\mid k\in\mathbb N\}\subseteq X_m$ and $x_k\to x$ in the topology of $X_m$.
\end{enumerate}
\end{theorem}

The argument for strict inductive limits comes from \cite[Chapter~13]{Treves1967}.

\begin{proof}
We will prove completeness using filters, since an $LF$ space need not be metrizable.

We will repeatedly use the following elementary observation. If $A$ is an absolutely convex set, $u,v\in A$ and $\alpha,\beta\in\mathbb K$ satisfy $|\alpha|+|\beta|\leq 1$, then $\alpha u+\beta v\in A$.

If $\alpha=0$ or $\beta=0$, the assertion follows from the balancedness of $A$. Now suppose that $\alpha\neq 0$ and $\beta\neq 0$, and define $s:=|\alpha|+|\beta|$. Since $0<s\leq 1$, we have
\[
\alpha u+\beta v
=
s\left(
\frac{|\alpha|}{s}\frac{\alpha}{|\alpha|}u
+
\frac{|\beta|}{s}\frac{\beta}{|\beta|}v
\right).
\]
Since $A$ is balanced, $\displaystyle\frac{\alpha}{|\alpha|}u,\displaystyle\frac{\beta}{|\beta|}v\in A$; since $\displaystyle\frac{|\alpha|}{s}+\displaystyle\frac{|\beta|}{s}=1$, the term in parentheses belongs to $A$ by convexity. Finally, since $s\leq 1$ and $A$ is balanced, we conclude that $\alpha u+\beta v\in A$.

\begin{enumerate}
 \item Let $x\in X\setminus\{0\}$. Since $X=\displaystyle\bigcup_{n=1}^{\infty}X_n$, there exists $m\in\mathbb N$ such that $x\in X_m$. Since $X_m$ is Fréchet, it is Hausdorff; therefore, there is a neighborhood of the origin in $X_m$ that does not contain $x$. By Proposition~\ref{prop: localmente convexo tiene una base de vecindades absolutamente convexas}, we can choose an absolutely convex neighborhood $U_m$ of the origin in $X_m$ such that $x\notin U_m$.

For $n<m$, define $U_n:=U_m\cap X_n$. Since the topology of $X_n$ agrees with the topology induced by $X_m$, each $U_n$ is an absolutely convex neighborhood of the origin in $X_n$. Moreover, for $2\leq n\leq m$ we have $U_n\cap X_{n-1}=U_{n-1}$.

We now construct $U_n$ inductively for $n>m$. Suppose that $U_{n-1}$ has been defined as an absolutely convex neighborhood of the origin in $X_{n-1}$. Since the topology of $X_{n-1}$ agrees with the topology induced by $X_n$, there exists an absolutely convex neighborhood $V_n$ of the origin in $X_n$ such that $V_n\cap X_{n-1}\subseteq U_{n-1}$. Define
\[
U_n:=
\operatorname{cobal}(V_n\cup U_{n-1})\footnotemark
=
\bigl\{
 \alpha v+\beta u
 \mid
 v\in V_n,\ u\in U_{n-1},\
 |\alpha|+|\beta|\leq 1
\bigr\}.
\]
By construction, $U_n$ is absolutely convex and contains $V_n$; in particular, it is a neighborhood of the origin in $X_n$.

We claim that $U_n\cap X_{n-1}=U_{n-1}$. The inclusion $U_{n-1}\subseteq U_n\cap X_{n-1}$ is immediate, since $U_{n-1}\subseteq U_n$ and $U_{n-1}\subseteq X_{n-1}$. To prove the other inclusion, take $z\in U_n\cap X_{n-1}$. Then there exist $v\in V_n$, $u\in U_{n-1}$, and scalars $\alpha,\beta\in\mathbb K$ such that $|\alpha|+|\beta|\leq 1$ and $z=\alpha v+\beta u$. If $\alpha=0$, then $z=\beta u\in U_{n-1}$, since $U_{n-1}$ is balanced. If $\alpha\neq0$, the preceding equality gives
\[
v=\alpha^{-1}(z-\beta u)\in X_{n-1},
\]
since $z,u\in X_{n-1}$. Therefore, $v\in V_n\cap X_{n-1}\subseteq U_{n-1}$. Since also $u\in U_{n-1}$ and $|\alpha|+|\beta|\leq1$, the preliminary observation about absolutely convex sets gives $z=\alpha v+\beta u\in U_{n-1}$. The equality is proved.

Now define
\[
U:=\displaystyle\bigcup_{n=1}^{\infty}U_n.
\]
The sequence $(U_n)_{n\in\mathbb N}$ is increasing, since $U_{n-1}=U_n\cap X_{n-1}\subseteq U_n$. It follows that $U$ is absolutely convex. If $y,z\in U$, there exist indices $r,s\in\mathbb N$ such that $y\in U_r$ and $z\in U_s$. Taking $N:=\displaystyle\max\{r,s\}$, we obtain $y,z\in U_N$; therefore, if $\alpha,\beta\in\mathbb K$ satisfy $|\alpha|+|\beta|\leq1$, then $\alpha y+\beta z\in U_N\subseteq U$, since $U_N$ is absolutely convex.

We show that $U\cap X_n=U_n$ for every $n\in\mathbb N$. If $r\geq n$, the equalities $U_j\cap X_{j-1}=U_{j-1}$ allow us to descend successively from $r$ to $n$:
\[
U_r\cap X_n
=
(U_r\cap X_{r-1})\cap X_n
=
U_{r-1}\cap X_n
=
\cdots
=
U_n.
\]
If, instead, $r<n$, then $U_r\subseteq U_n\subseteq X_n$, so $U_r\cap X_n=U_r\subseteq U_n$. Consequently,
\[
U\cap X_n
=
\left(\bigcup_{r=1}^{\infty}U_r\right)\cap X_n
=
\bigcup_{r=1}^{\infty}(U_r\cap X_n)
=
U_n.
\]
Thus, $U\cap X_n$ is a neighborhood of the origin in $X_n$ for every $n\in\mathbb N$. By definition of the topology $LF$, $U$ is a neighborhood of the origin in $X$. Applying the preceding equality to the index $m$ fixed at the outset, we obtain $U\cap X_m=U_m$. Since $x\in X_m$ and $x\notin U_m$, we conclude that $x\notin U$.

For each $x\in X\setminus\{0\}$, we have found a neighborhood of the origin that does not contain $x$. Therefore, $X$ is Hausdorff.
 \item First suppose that $B\subseteq X_m$ is topologically bounded in $X_m$. Let $U$ be a neighborhood of the origin in $X$. Since the inclusion $\iota_m\colon X_m\hookrightarrow X$ is continuous, $U\cap X_m$ is a neighborhood of the origin in $X_m$. By the topological boundedness of $B$ in $X_m$, there exists $\lambda>0$ such that
\[
B\subseteq \lambda(U\cap X_m)\subseteq \lambda U.
\]
Thus, $B$ is topologically bounded in $X$.

Conversely, suppose that $B$ is topologically bounded in $X$ and is not contained in any $X_m$. We recursively choose points $x_k\in B$ and indices $m_k\in\mathbb N$ as follows. Take $x_1\in B\setminus X_1$ and let $m_1$ be the least index such that $x_1\in X_{m_1}$. Once $m_{k-1}$ has been chosen, take $x_k\in B\setminus X_{m_{k-1}}$ and let $m_k$ be the least index such that $x_k\in X_{m_k}$. Then $(m_k)_{k\in\mathbb N}$ is strictly increasing and, for each $k$, $x_k\in X_{m_k}\setminus X_{m_k-1}$.

Since the topology of $X_{m_k-1}$ agrees with the topology induced by $X_{m_k}$ and $X_{m_k-1}$ is complete, being Fréchet, it follows that $X_{m_k-1}$ is a closed subspace of $X_{m_k}$. Since $x_k\notin X_{m_k-1}$, Corollary~\ref{cor: separacion de un subespacio y un abierto convexo, consecuencia de Hahn Banach} guarantees the existence of a continuous linear functional $f_k\in X_{m_k}'$ such that $f_k|_{X_{m_k-1}}=0$ and $f_k(x_k)=k$.

We construct a neighborhood $U$ of the origin in $X$ from an increasing sequence of absolutely convex neighborhoods $U_n\subseteq X_n$.

For $n=1$, fix an absolutely convex neighborhood $U_1$ of the origin in $X_1$. Suppose that $U_{n-1}$ has been constructed as an absolutely convex neighborhood of the origin in $X_{n-1}$. Since the topology of $X_{n-1}$ agrees with the topology induced by $X_n$, there exists an absolutely convex neighborhood $V_n$ of the origin in $X_n$ such that $V_n\cap X_{n-1}\subseteq U_{n-1}$. Define
\[
W_n:=
\operatorname{cobal}(V_n\cup U_{n-1})
=
\bigl\{
 \alpha v+\beta u
 \mid
 v\in V_n,\ u\in U_{n-1},\
 |\alpha|+|\beta|\leq 1
\bigr\}.
\]
Since $W_n$ contains $V_n$, it is a neighborhood of the origin in $X_n$; moreover, by construction, it is absolutely convex. The same argument used in the preceding part shows that $W_n\cap X_{n-1}=U_{n-1}$.

Now define
\[
U_n:=
\begin{cases}
W_n, & \text{if } n\neq m_k \text{ for every }k,\\[4pt]
W_n\cap\{x\in X_n\mid |f_k(x)|\leq 1\},
& \text{if } n=m_k.
\end{cases}
\]
Each $U_n$ is absolutely convex. Moreover, if $n=m_k$, the set $\{x\in X_n\mid |f_k(x)|\leq 1\}$ is a neighborhood of the origin in $X_n$, since it contains the open set $\{x\in X_n\mid |f_k(x)|<1\}$. Therefore, each $U_n$ is an absolutely convex neighborhood of the origin in $X_n$.

The equality $U_n\cap X_{n-1}=U_{n-1}$ is preserved at every step. If $n\neq m_k$ for every $k$, we simply have $U_n\cap X_{n-1}=W_n\cap X_{n-1}=U_{n-1}$. If $n=m_k$, then
\[
\begin{aligned}
U_{m_k}\cap X_{m_k-1}
&=
W_{m_k}\cap\{x\in X_{m_k}\mid |f_k(x)|\leq 1\}\cap X_{m_k-1}\\
&=
W_{m_k}\cap X_{m_k-1}
=
U_{m_k-1},
\end{aligned}
\]
since $f_k$ vanishes on $X_{m_k-1}$.

Since $U_{n-1}=U_n\cap X_{n-1}\subseteq U_n$, the sequence $(U_n)_{n\in\mathbb N}$ is increasing. Define $U:=\displaystyle\bigcup_{n=1}^{\infty}U_n$. As the union of an increasing sequence of absolutely convex sets, $U$ is absolutely convex; moreover, the preceding equalities give $U\cap X_n=U_n$ for every $n\in\mathbb N$. By definition of the topology $LF$, $U$ is a neighborhood of the origin in $X$.

Since $B$ is topologically bounded in $X$, there exists $\lambda>0$ such that $B\subseteq \lambda U$. In particular, $x_k\in\lambda U$ for every $k$. Since $x_k\in X_{m_k}$, we have
\[
x_k\in \lambda U\cap X_{m_k}
=
\lambda(U\cap X_{m_k})
=
\lambda U_{m_k}.
\]
By definition of $U_{m_k}$, every element $y\in U_{m_k}$ satisfies $|f_k(y)|\leq 1$. Since $\displaystyle\frac{x_k}{\lambda}\in U_{m_k}$,
\[
\left|f_k\!\left(\frac{x_k}{\lambda}\right)\right|\leq 1,
\]
and therefore $|f_k(x_k)|\leq\lambda$. However, by construction, $f_k(x_k)=k$. Thus, $k\leq\lambda$ for every $k\in\mathbb N$, which is impossible.

This contradiction shows that $B$ must be contained in some $X_m$. It remains to show that it is topologically bounded in $X_m$. Let $V$ be a neighborhood of the origin in $X_m$. Since the topology of $X_m$ agrees with the topology induced by $X$, there exists a neighborhood $W$ of the origin in $X$ such that $V=W\cap X_m$. Since $B$ is topologically bounded in $X$, there exists $\lambda>0$ such that $B\subseteq \lambda W$. Since also $B\subseteq X_m$,
\[
B\subseteq \lambda W\cap X_m
=
\lambda(W\cap X_m)
=
\lambda V.
\]
Therefore, $B$ is topologically bounded in $X_m$.

 \item We prove completeness using filters, since $X$ is not necessarily metrizable. Recall that a filter $\mathscr F$ on $X$ is a family of nonempty subsets closed under finite intersections and passage to supersets. We say that $\mathscr F$ is Cauchy if, for every neighborhood $U$ of the origin, there exists $M\in\mathscr F$ such that $M-M\subseteq U$; the space is complete if every Cauchy filter converges.

 Let $\mathscr F$ be a Cauchy filter on $X$. The sets
 \[
 M+V,
 \qquad
 M\in\mathscr F,
 \quad
 V\text{ a neighborhood of the origin in }X,
 \]
 form a filter base: the intersection of two of them contains $(M_1\cap M_2)+(V_1\cap V_2)$. Denote the filter they generate by $\mathscr G$. Each basic set $M+V$ belongs to $\mathscr F$, since it contains $M$; therefore, $\mathscr F$ is finer than $\mathscr G$. Moreover, $\mathscr G$ is Cauchy. Indeed, given a neighborhood $U$ of the origin, choose an absolutely convex neighborhood $V$ such that $V+V-V\subseteq U$ and a set $M\in\mathscr F$ such that $M-M\subseteq V$. Then
 \[
 (M+V)-(M+V)
 \subseteq
 V+V-V
 \subseteq U.
 \]

 We claim that there exists $p\in\mathbb N$ such that
 \[
 A\cap X_p\neq\varnothing
 \qquad
 \text{for every }A\in\mathscr G.
 \tag{\(*\)}
 \]
 Suppose otherwise. For each $n$ we could choose a basic set $A_n=M_n+V_n\in\mathscr G$ such that $A_n\cap X_n=\varnothing$, with $M_n\in\mathscr F$. By successively shrinking the neighborhoods, we may assume that each $V_n$ is absolutely convex and that $V_{n+1}\subseteq V_n$. Define
 \[
 W_n
 :=
 \operatorname{conv}
 \left(
 V_n\cup
 \bigcup_{k<n}(V_k\cap X_k)
 \right).
 \]
 We have $(M_n+W_n)\cap X_n=\varnothing$. For $n=1$ this follows from $W_1=V_1$. If $n\geq2$, every element of $W_n$ can be written as $ty+(1-t)z$, with $0\leq t\leq1$, $y\in V_n$, and $z\in X_{n-1}$, after replacing $z$ by an element of the convex hull of $\displaystyle\bigcup_{k<n}(V_k\cap X_k)$. If $m+ty+(1-t)z\in X_n$, with $m\in M_n$, then $m+ty\in(M_n+V_n)\cap X_n$, because $V_n$ is balanced, contradicting the choice of $A_n$.

 Now set
 \[
 W
 :=
 \operatorname{conv}
 \left(
 \bigcup_{k=1}^{\infty}(V_k\cap X_k)
 \right).
 \]
 The set $W$ is absolutely convex and $W\cap X_k$ contains the neighborhood $V_k\cap X_k$ of the origin in $X_k$; by definition of the topology $LF$, $W$ is a neighborhood of the origin in $X$. The monotonicity of $(V_k)$ also shows that $W\subseteq W_n$ for every $n$: the terms with index less than $n$ appear in the definition of $W_n$, whereas those with index greater than or equal to $n$ are contained in $V_n$.

 Since $\mathscr G$ is Cauchy, there exists $B\in\mathscr G$ such that $B-B\subseteq W$. The set $B$ also belongs to $\mathscr F$, and therefore $B\cap M_n\neq\varnothing$ for every $n$. If $b_n\in B\cap M_n$, then
 \[
 B
 \subseteq
 b_n+(B-B)
 \subseteq
 M_n+W_n.
 \]
 Consequently, $B\cap X_n=\varnothing$ for every $n$, which is impossible because $X=\displaystyle\bigcup_nX_n$ and a filter does not contain the empty set. This proves assertion~\textup{\((*)\)}.

 Fix $p$ as in~\textup{\((*)\)}. The sets $A\cap X_p$, with $A\in\mathscr G$, form a base for a filter $\mathscr G_p$ on $X_p$. This filter is Cauchy: if $U_p$ is a neighborhood of the origin in $X_p$, equality of the original topology of $X_p$ and the one induced by $X$ allows us to choose a neighborhood $U$ of the origin in $X$ such that $U\cap X_p\subseteq U_p$, and it suffices to apply the Cauchy property of $\mathscr G$. Since $X_p$ is Fréchet, there exists $x\in X_p$ such that $\mathscr G_p$ converges to $x$.

 This implies that $x$ is a cluster point of $\mathscr G$. A Cauchy filter with a cluster point converges to it: given a neighborhood $U$ of the origin, choose a balanced neighborhood $V$ with $V+V\subseteq U$ and $A\in\mathscr G$ with $A-A\subseteq V$; since $A\cap(x+V)\neq\varnothing$, if $a$ belongs to this intersection, then $A\subseteq x+V+V\subseteq x+U$. Thus, $\mathscr G$ converges to $x$. Since $\mathscr F$ is finer than $\mathscr G$, $\mathscr F$ also converges to $x$. This proves that $X$ is complete.

 \item First suppose that $x_k\to x$ in $X$. Then $(x_k)_{k\in\mathbb N}$ is a Cauchy sequence, and by Proposition~\ref{prop: sucesion de Cauchy tiene rango topologicamente acotado}, the set $\{x_k\mid k\in\mathbb N\}$ is topologically bounded in $X$. Since the singleton $\{x\}$ is also topologically bounded, the union $B:=\{x\}\cup\{x_k\mid k\in\mathbb N\}$ is topologically bounded in $X$. By part (b), there exists $m\in\mathbb N$ such that $B\subseteq X_m$. Since the topology of $X_m$ agrees with the topology induced by $X$, it follows that $x_k\to x$ in $X_m$.

 Conversely, if there exists $m\in\mathbb N$ such that $\{x\}\cup\{x_k\mid k\in\mathbb N\}\subseteq X_m$ and $x_k\to x$ in $X_m$, then $x_k\to x$ in $X$ by continuity of the canonical inclusion $X_m\hookrightarrow X$.
\end{enumerate}
\end{proof}

\footnotetext{
If $A\subseteq X$ is a subset of a vector space, its absolutely convex hull is defined by
\[
\operatorname{cobal}(A)
=
\left\{
\displaystyle\sum_{j=1}^{r}\lambda_j x_j
\;\Biggr|\;
x_j\in A,\
\displaystyle\sum_{j=1}^{r}|\lambda_j|\le 1
\right\}.
\]
It is the smallest absolutely convex subset of $X$ containing $A$.
}

\section{Distributions: the dual of \texorpdfstring{$\mathcal{D}(\Omega)$}{D(Omega)}}
The concept of a distribution will be obtained from the dual of $\mathcal{D}(\Omega)$. To this end, we first study some properties of continuous linear functionals defined on $\mathcal{D}(\Omega)$.
\begin{theorem}\label{teo: equivalencias de continuidad en D'(Omega)}\index{equivalent conditions for continuity in D'(Omega)}
 Let $\Omega\subseteq \mathbb{R}^{n}$ be an open set and let $T\colon \mathcal{D}(\Omega)\longrightarrow \mathbb{R}$ be linear. Then the following properties are equivalent:
 \begin{enumerate}
 \item $T$ is continuous.
 \item $T$ is bounded.
 \item If $\{\phi_{\ell}\}_{\ell\in \mathbb{N}}$ converges to $\phi\in \mathcal{D}(\Omega)$ with respect to $\tau$, then
 \[\lim_{\ell\to \infty}T(\phi_{\ell})=T(\phi).\]
 \item $T\restriction_{\mathcal{D}_{K}(\Omega)}$ is continuous for every compact set $K\subseteq \Omega$.
 \item For every compact set $K\subseteq \Omega$, there exist $j\in \mathbb{N}_{0}$ and $c_{K,j}>0$ such that
 \begin{equation}\label{eq: continuidad T distribucion constante compacto}
 |T(\phi)|\leq c_{K,j}\|\phi\|_{C^{j}(K)}\qquad \forall \phi\in \mathcal{D}_{K}(\Omega). \end{equation}
 \end{enumerate}
\end{theorem}

\begin{proof}

We show that (a) implies (b). If $T$ is continuous, then $T$ is bounded by Theorem~\ref{teo: relacion continuidad y acotacion en espacios vectoriales topologicos}.

For (b)$\Rightarrow$(c), suppose that $T$ is bounded and let $\{\phi_{\ell}\}_{\ell\in\mathbb N}$ be a sequence in $\mathcal D(\Omega)$ such that $\phi_{\ell}\to\phi$ with respect to $\tau$. By part (d) of Theorem~\ref{teo: propiedades LF}, there exists a compact set $K\subseteq\Omega$ such that
\[
 \{\phi\}\cup\{\phi_{\ell}\mid \ell\in\mathbb N\}
 \subseteq\mathcal D_K(\Omega)
\]
and $\phi_{\ell}\to\phi$ in the Fréchet topology of $\mathcal D_K(\Omega)$.

The restriction
\[
 T\restriction_{\mathcal D_K(\Omega)}
 \colon\mathcal D_K(\Omega)\longrightarrow\mathbb R
\]
is bounded. Indeed, the inclusion $\mathcal D_K(\Omega)\hookrightarrow\mathcal D(\Omega)$ is continuous and therefore maps topologically bounded sets to topologically bounded sets; one then applies the boundedness of $T$. Since $\mathcal D_K(\Omega)$ is metrizable, Theorem~\ref{teo: relacion continuidad y acotacion en espacios vectoriales topologicos} shows that this restriction is continuous. Consequently,
\[
 T(\phi_{\ell})\longrightarrow T(\phi).
\]

(c)$\Rightarrow$(d). Let $K\subseteq \Omega$ be compact. Suppose that $\{\phi_{\ell}\}_{\ell\in \mathbb{N}}\subseteq \mathcal{D}_{K}(\Omega)$ and that
\[
d_{K}(\phi_{\ell},\phi)\to 0.
\]

By Proposition~\ref{prop: convergencia y Cauchy por seminormas}, this means that $\phi_{\ell}\to \phi$ in $(\mathcal{D}_{K}(\Omega),\tau_{K})$. Since the inclusion
\[
\mathcal{D}_{K}(\Omega)\hookrightarrow \mathcal{D}(\Omega)
\]
is continuous, it follows that $\phi_{\ell}\to \phi$ with respect to $\tau$. By (c),
\[
T(\phi_{\ell})\to T(\phi).
\]

Thus, $T\restriction_{\mathcal{D}_{K}(\Omega)}$ is sequentially continuous. Since $\mathcal{D}_{K}(\Omega)$ is metrizable, we conclude that $T\restriction_{\mathcal{D}_{K}(\Omega)}$ is continuous.

(d)$\Rightarrow$(a). Let $V\subseteq\mathbb R$ be an absolutely convex neighborhood of the origin and set
\[
 U:=T^{-1}(V).
\]
The set $U$ is absolutely convex. For every compact set $K\subseteq\Omega$, continuity of $T\restriction_{\mathcal D_K(\Omega)}$ implies that
\[
 U\cap\mathcal D_K(\Omega)
 =\bigl(T\restriction_{\mathcal D_K(\Omega)}\bigr)^{-1}(V)
\]
is a neighborhood of the origin in $\mathcal D_K(\Omega)$. By the description of neighborhoods of the origin for the locally convex final topology in Theorem~\ref{teo: topologia-final-lc}, $U$ is a neighborhood of the origin in $\mathcal D(\Omega)$. Thus, $T$ is continuous at the origin and, by linearity, at every point.

(d)$\Leftrightarrow$(e). Suppose that (d) holds and let $K\subseteq \Omega$ be compact. Since $T\restriction_{\mathcal{D}_{K}(\Omega)}$ is continuous at the origin, given $\varepsilon=1$, there exist $j\in \mathbb{N}_{0}$ and $n\in \mathbb{N}$ such that
\[
V_{K,j,n}\subseteq T^{-1}((-1,1)).
\]

This means that
\[
|T(\phi)|\leq 1\qquad \forall \phi\in \mathcal{D}_{K}(\Omega)\ \text{with}\ \|\phi\|_{C^{j}(K)}<\frac{1}{n}.
\]

If $\phi\in \mathcal{D}_{K}(\Omega)$ and $\phi\neq 0$, then
\[
\left\|\frac{1}{2n}\frac{\phi}{\|\phi\|_{C^{j}(K)}}\right\|_{C^{j}(K)}<\frac{1}{n}.
\]

Therefore,
\[
\left|T\left(\frac{1}{2n}\frac{\phi}{\|\phi\|_{C^{j}(K)}}\right)\right|\leq 1.
\]

By linearity of $T$, it follows that
\[
|T(\phi)|\leq 2n\|\phi\|_{C^{j}(K)}.
\]

This proves (e). Conversely, suppose that (e) holds and let $K\subseteq \Omega$ be compact. We wish to prove that
\[
T\restriction_{\mathcal{D}_{K}(\Omega)}\colon \mathcal{D}_{K}(\Omega)\longrightarrow \mathbb{R}
\]
is continuous.

By hypothesis, there exist $j\in\mathbb{N}_{0}$ and $c_{K}>0$ such that
\[
|T(\phi)|\leq c_{K,j}\|\phi\|_{C^{j}(K)}\qquad \forall \phi\in \mathcal{D}_{K}(\Omega).
\]

It suffices to show that $T\restriction_{\mathcal{D}_{K}(\Omega)}$ is continuous at $0$. Since $T$ is linear, continuity at the origin implies continuity at every point.

Let $\varepsilon>0$. By the Archimedean property, there exists $n\in\mathbb N$ such that $n>\displaystyle\frac{c_{K,j}}{\varepsilon}$, and therefore
\[
\frac{c_{K,j}}{n}<\varepsilon.
\]

Now consider $\phi\in V_{K,j,n}$. By definition of $V_{K,j,n}$, this means that
\[
\|\phi\|_{C^{j}(K)}<\frac{1}{n}.
\]

Applying the bound given by (e), we obtain
\[
|T(\phi)|\leq c_{K,j}\|\phi\|_{C^{j}(K)}<\frac{c_{K,j}}{n}<\varepsilon.
\]

Therefore,
\[
\phi\in V_{K,j,n}\quad \Longrightarrow \quad T(\phi)\in (-\varepsilon,\varepsilon),
\]
and consequently
\[
V_{K,j,n}\subseteq \left(T\restriction_{\mathcal{D}_{K}(\Omega)}\right)^{-1}((-\varepsilon,\varepsilon)).
\]

Since $V_{K,j,n}$ is a neighborhood of $0$ in $\mathcal{D}_{K}(\Omega)$, we conclude that
\[
\left(T\restriction_{\mathcal{D}_{K}(\Omega)}\right)^{-1}((-\varepsilon,\varepsilon))
\]
is a neighborhood of $0$. This proves that $T\restriction_{\mathcal{D}_{K}(\Omega)}$ is continuous at $0$ and, by linearity, is continuous throughout $\mathcal{D}_{K}(\Omega)$.
\end{proof}
\begin{definition}\label{def:distribuciones-y-fourier-el-espacio-de-distribuciones-definidas-sobre}\index{distribution!on a Euclidean open set}\index{space of distributions}
 Let $\Omega\subseteq \mathbb{R}^{n}$ be open. The continuous dual of $(\mathcal{D}(\Omega),\tau)$, denoted by $\mathcal{D}'(\Omega)$, is called \textit{\textbf{the space of distributions defined on $\Omega$}}, and its elements are called \textit{\textbf{distributions}}. We will frequently use the duality notation $\langle T,\phi\rangle$ to denote $T(\phi)$. The space $\mathcal{D}'(\Omega)$ is equipped with the weak-$*$ topology. Thus, a sequence $\{T_{k}\}_{k\in \mathbb{N}}\subseteq \mathcal{D}'(\Omega)$ converges to $T\in \mathcal{D}'(\Omega)$ if \[T_{k}(\phi)\to T(\phi),\qquad \forall \phi\in \mathcal{D}(\Omega).\]

 In this case, we say that $\{T_{k}\}_{k\in \mathbb{N}}$ converges to $T$ \textit{in the sense of distributions} or \textit{in the distributional sense}.
 \end{definition}
\begin{example}\label{ej:distribuciones-y-fourier-funcional-definido-distribucion-primero-observamos-lineal}
 If $\psi\in C^{\infty}(\Omega)$ and $T\in\mathcal{D}'(\Omega)$, then the functional $\psi T\colon \mathcal{D}(\Omega)\longrightarrow \mathbb{R}$ defined by
 \[
 (\psi T)(\phi):=T(\psi\phi),\qquad \phi\in\mathcal{D}(\Omega),
 \]
 is a distribution.

 First observe that $\psi T$ is linear. If $\phi_{1},\phi_{2}\in\mathcal{D}(\Omega)$ and $a,b\in\mathbb{R}$, then
 \[
 (\psi T)(a\phi_{1}+b\phi_{2})
 =
 T\big(\psi(a\phi_{1}+b\phi_{2})\big)
 =
 T(a\psi\phi_{1}+b\psi\phi_{2})
 =
 a(\psi T)(\phi_{1})+b(\psi T)(\phi_{2}).
 \]

 We now show that $\psi T$ is continuous. Let $K\subseteq\Omega$ be compact. Since $T$ is a distribution, by part (e) of Theorem~\ref{teo: equivalencias de continuidad en D'(Omega)}, there exist $j\in\mathbb{N}_{0}$ and $c_{K}>0$ such that
 \[
 |T(\eta)|\leq c_{K}\|\eta\|_{C^{j}(K)},
 \qquad
 \forall \eta\in\mathcal D_{K}(\Omega).
 \]

 Now let $\phi\in\mathcal D_{K}(\Omega)$. Then $\psi\phi\in\mathcal D_{K}(\Omega)$ and therefore
 \[
 |(\psi T)(\phi)|
 =
 |T(\psi\phi)|
 \leq
 c_{K}\|\psi\phi\|_{C^{j}(K)}.
 \]

 To estimate this seminorm, we use the Leibniz rule for partial derivatives with multi-indices (Proposition~\ref{prop: propiedades parciales multiindices incluye leibniz multiindices}). If $|\alpha|\leq j$, then
 \[
 \partial^{\alpha}(\psi\phi)
 =
 \displaystyle\sum_{\beta\leq\alpha}
 \binom{\alpha}{\beta}
 (\partial^{\beta}\psi)
 (\partial^{\alpha-\beta}\phi).
 \]

 Taking the supremum over $K$ gives
 \[
 \|\partial^{\alpha}(\psi\phi)\|_{C^0(K)}
 \leq
 \displaystyle\sum_{\beta\leq\alpha}
 \binom{\alpha}{\beta}
 \|\partial^{\beta}\psi\|_{C^0(K)}
 \|\partial^{\alpha-\beta}\phi\|_{C^0(K)}.
 \]

 Since $\psi\in C^{\infty}(\Omega)$ and $K$ is compact, each derivative $\partial^{\beta}\psi$ is bounded on $K$. Define
 \[
 M_{\beta}:=\|\partial^{\beta}\psi\|_{C^0(K)}.
 \]

 Then
 \[
 \|\partial^{\alpha}(\psi\phi)\|_{C^0(K)}
 \leq
 \displaystyle\sum_{\beta\leq\alpha}
 \binom{\alpha}{\beta}
 M_{\beta}
 \|\partial^{\alpha-\beta}\phi\|_{C^0(K)}.
 \]

 Summing over $|\alpha|\leq j$ shows that there exists a constant $C_{K}>0$ such that
 \[
 \|\psi\phi\|_{C^{j}(K)}\leq C_{K}\|\phi\|_{C^{j}(K)}.
 \]

 Substituting into the preceding estimate yields
 \[
 |(\psi T)(\phi)|
 \leq
 c_{K}C_{K}\|\phi\|_{C^{j}(K)}.
 \]

 Applying part (e) of Theorem~\ref{teo: equivalencias de continuidad en D'(Omega)} again, we conclude that $\psi T$ is a distribution, that is,
 \[
 \psi T\in\mathcal D'(\Omega).
 \]
\end{example}
\begin{example}\label{ej:distribuciones-y-fourier-funcional-distribucion-efecto-lineal-linealidad-integral}
Let $f\in L_{\text{loc}}^{1}(\Omega)$. Define the functional
\[
T_{f}\colon \mathcal{D}(\Omega)\longrightarrow \mathbb{R},\qquad
T_{f}(\phi):=\int_{\Omega} f(x)\phi(x)\,dx.
\]
Then $T_f$ is a distribution. Linearity follows from linearity of the integral. To prove continuity, let $K\subseteq \Omega$ be compact and let $\phi\in\mathcal{D}_{K}(\Omega)$. Then
\[
|T_{f}(\phi)|
=
\left|\int_{K} f(x)\phi(x)\,dx\right|
\leq
\int_{K} |f(x)||\phi(x)|\,dx.
\]

Since $|\phi(x)|\leq \|\phi\|_{C^0(K)}$ for every $x\in K$, we obtain
\[
|T_{f}(\phi)|
\leq
\|\phi\|_{C^0(K)}\int_{K}|f(x)|\,dx.
\]

Defining
\[
c_{K}:=\int_{K}|f(x)|\,dx,
\]
we have
\[
|T_{f}(\phi)|\leq c_{K}\|\phi\|_{C^0(K)}.
\]

By part (e) of Theorem~\ref{teo: equivalencias de continuidad en D'(Omega)}, we conclude that $T_{f}$ is a distribution.
\end{example}
\begin{example}\label{ej:distribuciones-y-fourier-delta-de-dirac}
Let $x_{0}\in\Omega$. The functional
\[
\delta_{x_{0}}\colon \mathcal{D}(\Omega)\longrightarrow \mathbb{R},\qquad
\delta_{x_{0}}(\phi):=\phi(x_{0}),
\]
is called the \textbf{Dirac delta} centered at $x_{0}$. Linearity is immediate. We show that it is continuous. Let $K\subseteq\Omega$ be compact with $x_{0}\in K$ and let $\phi\in\mathcal{D}_{K}(\Omega)$. Then
\[
|\delta_{x_{0}}(\phi)|
=
|\phi(x_{0})|
\leq
\|\phi\|_{C^0(K)}.
\]

Therefore,
\[
|\delta_{x_{0}}(\phi)|\leq \|\phi\|_{C^0(K)}.
\]

If $x_0\notin K$, every function in $\mathcal D_K(\Omega)$ vanishes at $x_0$, and the same bound holds because the functional is zero.

Applying part (e) of Theorem~\ref{teo: equivalencias de continuidad en D'(Omega)}, we conclude that $\delta_{x_{0}}$ is a distribution.
\end{example}
An important operation on distributions is the \textit{restriction} of a distribution: \begin{definition}\label{def:distribuciones-y-fourier-r-n-abierto-v-d-r}\index{restriction}
 Let \(\Omega\subseteq \mathbb{R}^{n}\) be an open set and let \(v\in \mathcal{D}'(\mathbb{R}^{n})\). We define the restriction of \(v\) to \(\Omega\) to be the distribution \(v\restriction_{\Omega}\in \mathcal{D}'(\Omega)\) given by
 \[
 \langle v\restriction_{\Omega},\varphi\rangle
 :=
 \langle v,\tilde{\varphi}\rangle,
 \qquad \forall \varphi\in \mathcal{D}(\Omega),
 \]
 where \(\tilde{\varphi}\in \mathcal{D}(\mathbb{R}^{n})\) is the extension of \(\varphi\) by zero outside \(\Omega\).
\end{definition}
\section{Order of a distribution}
In this section we define the order of a distribution. The preceding examples show that the integer $j\in \mathbb{N}_{0}$ in \eqref{eq: continuidad T distribucion constante compacto} may change with the choice of distribution and compact set. If the same integer works for every compact set $K\subseteq \Omega$, then the smallest integer for which \eqref{eq: continuidad T distribucion constante compacto} holds is known as the \textit{order of the distribution} $T$.
\begin{definition}\label{def:distribuciones-y-fourier-norma-distribucion-orden-menor-igual}\index{distribution!order}
Let $T\in\mathcal{D}'(\Omega)$ and fix $m\in\mathbb N_0$. We say that $T$ has \textit{order less than or equal to $m$} if for every compact set $K\subseteq\Omega$ there exists a constant $c_{K}>0$ with the property that
\[
|T(\phi)|\leq c_{K}\|\phi\|_{C^m(K)}
\qquad
\forall \phi\in\mathcal{D}_{K}(\Omega).
\]

If such an integer $m$ exists, the smallest one is called the \textit{order of the distribution} $T$ and is denoted by $\operatorname{ord}(T)$. If no such integer exists, we say that $T$ is a distribution of \textit{infinite order}.
\end{definition}
\begin{example}\label{ej:distribuciones-y-fourier-ejemplos-anteriores-permiten-identificar-orden-varias}
The preceding examples allow us to identify the order of several important distributions.

\begin{enumerate}

\item If $f\in L^{1}_{\mathrm{loc}}(\Omega)$, then the distribution
\[
T_f(\phi)=\int_\Omega f(x)\phi(x)\,dx
\]
has order $0$.

If $K\subseteq\Omega$ is compact and $\phi\in\mathcal{D}_K(\Omega)$, then
\[
|T_f(\phi)|
=
\left|\int_K f(x)\phi(x)\,dx\right|
\le
\int_K |f(x)||\phi(x)|\,dx
\le
\|\phi\|_{C^0(K)}\int_K |f(x)|\,dx.
\]

Defining
\[
c_K:=\int_K |f(x)|\,dx,
\]
we obtain
\[
|T_f(\phi)|\le c_K\|\phi\|_{C^0(K)}.
\]

Therefore, the integer $m=0$ works in the preceding definition, and we conclude that
\[
\operatorname{ord}(T_f)=0.
\]

\item For $x_0\in\Omega$, the Dirac delta
\[
\delta_{x_0}(\phi)=\phi(x_0)
\]
is also a distribution of order $0$.

If $K\subseteq\Omega$ is compact and $x_0\in K$, then for $\phi\in\mathcal{D}_K(\Omega)$
\[
|\delta_{x_0}(\phi)|
=
|\phi(x_0)|
\le
\|\phi\|_{C^0(K)}.
\]

Therefore,
\[
|\delta_{x_0}(\phi)|\le \|\phi\|_{C^0(K)},
\]
and, if $x_0\notin K$, the functional is zero. Thus the integer $m=0$ satisfies the condition in the definition, implying
\[
\operatorname{ord}(\delta_{x_0})=0.
\]

\item If $\psi\in C^\infty(\Omega)$ and $T\in\mathcal{D}'(\Omega)$ has order $m$, then $\psi T$ also has order less than or equal to $m$.

The estimate in the preceding example provides a constant $C_K>0$ such that
\[
\|\psi\phi\|_{C^m(K)}\le C_K\|\phi\|_{C^m(K)}.
\]

Since $T$ has order $m$, there exists $c_K>0$ such that
\[
|T(\eta)|\le c_K\|\eta\|_{C^m(K)}.
\]

Taking $\eta=\psi\phi$, we obtain
\[
|(\psi T)(\phi)|
=
|T(\psi\phi)|
\le
c_K\|\psi\phi\|_{C^m(K)}
\le
c_KC_K\|\phi\|_{C^m(K)}.
\]

Therefore, $\psi T$ satisfies the definition with the same integer $m$, and thus
\[
\operatorname{ord}(\psi T)\le m.
\]

\item Now let $\mu$ be a signed Radon measure\footnotemark[1] on $\Omega$. Define
\[
T_\mu(\phi)=\int_\Omega \phi\,d\mu,
\qquad
\phi\in\mathcal{D}(\Omega).
\]

First observe that $T_\mu$ is linear by linearity of the integral.

We show that it is continuous. Let $K\subseteq\Omega$ be compact and let $\phi\in\mathcal{D}_K(\Omega)$. Since $\phi$ vanishes outside $K$, we have
\[
T_\mu(\phi)=\int_K \phi\,d\mu.
\]

Let $|\mu|$ be the total variation\footnotemark[2] of $\mu$. Then
\[
|T_\mu(\phi)|
=
\left|\int_K \phi\,d\mu\right|
\le
\int_K |\phi|\,d|\mu|.
\]

Since $|\phi(x)|\le\|\phi\|_{C^0(K)}$ for every $x\in K$, we obtain
\[
|T_\mu(\phi)|
\le
\|\phi\|_{C^0(K)}\int_K d|\mu|
=
|\mu|(K)\,\|\phi\|_{C^0(K)}.
\]

Defining
\[
c_K:=|\mu|(K),
\]
we obtain
\[
|T_\mu(\phi)|\le c_K\|\phi\|_{C^0(K)}.
\]

By the definition of the order of a distribution, we conclude that
\[
\operatorname{ord}(T_\mu)=0.
\]

\end{enumerate}
\end{example}

\footnotetext[1]{
A \textit{signed measure} on $\Omega$ is a function $\mu\colon \mathcal{B}(\Omega)\longrightarrow[-\infty,\infty]$ defined on the Borel $\sigma$-algebra $\mathcal{B}(\Omega)$ such that $\mu(\emptyset)=0$, taking at most one of the values $\infty$ and $-\infty$, and satisfying $\mu\Bigl(\bigcup_{k=1}^{\infty}E_k\Bigr)
=
\displaystyle\sum_{k=1}^{\infty}\mu(E_k)$ for every countable family $\{E_k\}_{k\in\mathbb{N}}\subseteq\mathcal{B}(\Omega)$ of disjoint sets. This is the global notion of Definition~\ref{def:b4-espacios-medida-medida-conjunto-no-vacio-algebra}. To be Radon, its positive and negative variations must additionally be Radon; see Definition~\ref{def:b4-espacios-medida-radon-con-signo}.}

\footnotetext[2]{
Let $\mu$ be a signed measure. Its \textit{total variation} is the positive measure $|\mu|\colon \mathcal{B}(\Omega)\longrightarrow[0,\infty]$ defined by
\[
|\mu|(E)
=
\sup
\Biggl\{
\displaystyle\sum_{k=1}^{\infty} |\mu(E_k)|
\;\Biggr|\;
E=\bigcup_{k=1}^{\infty}E_k,\;
E_k\in\mathcal{B}(\Omega),\;
E_i\cap E_j=\emptyset,\,i\neq j
\Biggr\}.
\]
Equivalently, if $\mu=\mu^{+}-\mu^{-}$ is the Jordan decomposition of $\mu$, then $|\mu|=\mu^{+}+\mu^{-}$.}
It turns out that all distributions of order $0$ can be identified with measures.
\begin{theorem}\label{teo: distribuciones orden cero y medidas}\index{distributions of order zero and measures}
 Let $\Omega\subseteq \mathbb{R}^{n}$ be an open set and let $T\in \mathcal{D}'(\Omega)$.
 \begin{enumerate}
 \item If $T$ is positive, that is, if
 \[
 T(\phi)\geq 0,\qquad \forall \phi\in \mathcal{D}(\Omega)\text{ with }\phi\geq 0,
 \]
 then there exists a unique positive Radon measure $\mu\colon \mathcal{B}(\Omega)\longrightarrow [0,\infty]$ such that
 \[
 T(\phi)=\int_{\Omega}\phi\,d\mu,\qquad \forall \phi\in \mathcal{D}(\Omega).
 \]

 \item If $T$ has order $0$, then there exist two positive Radon measures
 \[
 \mu_{1},\mu_{2}\colon \mathcal{B}(\Omega)\longrightarrow [0,\infty]
 \]
 such that
 \[
 T(\phi)=\int_{\Omega}\phi\,d\mu_{1}-\int_{\Omega}\phi\,d\mu_{2},\qquad \forall \phi\in \mathcal{D}(\Omega).
 \]
 \end{enumerate}
\end{theorem}

\begin{proof}
Recall that $C_c(\Omega)$ is equipped with the strict inductive limit topology
\[
C_c(\Omega)=\varinjlim_{i\in\mathbb N} C_{K_i}(\Omega),
\]
where $\{K_i\}_{i\in\mathbb N}$ is an increasing exhaustion of $\Omega$ by compact sets and
\[
C_{K_i}(\Omega):=\{u\in C_c(\Omega)\mid \supp(u)\subseteq K_i\}
\]
is equipped with the uniform norm
\[
\|u\|_{K_i,\infty}:=\|u\|_{C^0(K_i)}.
\]
In particular, a linear functional $L\colon C_c(\Omega)\longrightarrow\mathbb R$ is continuous if and only if for every compact set $K\subset\Omega$ there exists a constant $C_K>0$ such that
\[
|L(u)|\le C_K\|u\|_{K,\infty},
\qquad \forall u\in C_c(\Omega)\text{ with }\supp(u)\subset K .
\]

\begin{enumerate}

\item First suppose that $T$ is positive. We will show that $T$ has order $0$.

Let $K\subset\Omega$ be compact. Take a bounded open set $U\subset\Omega$ such that
\[
K\subset U\subset\overline U\subset\Omega .
\]
Consider a cutoff function $\psi\in C^\infty(\Omega)$ such that
\[
\psi\equiv1 \text{ in }K, \qquad
\supp(\psi)\subset U, \qquad
0\le\psi\le1 .
\]
In particular, $\psi\in\mathcal D(\Omega)$.

Now let $\phi\in\mathcal D_K(\Omega)$. Since $\psi\equiv1$ on $K$ and $\supp(\phi)\subset K$, we have
\[
|\phi(x)|\le\|\phi\|_{C^0(K)}\psi(x),
\qquad \forall x\in\Omega .
\]
It follows that
\[
\|\phi\|_{C^0(K)}\psi-\phi\ge0,
\qquad
\|\phi\|_{C^0(K)}\psi+\phi\ge0 .
\]
Since $T$ is positive,
\[
T(\|\phi\|_{C^0(K)}\psi-\phi)\ge0,
\qquad
T(\|\phi\|_{C^0(K)}\psi+\phi)\ge0 .
\]
By linearity,
\[
-\|\phi\|_{C^0(K)}T(\psi)\le T(\phi)\le \|\phi\|_{C^0(K)}T(\psi),
\]
and therefore
\[
|T(\phi)|\le T(\psi)\|\phi\|_{C^0(K)},
\qquad \forall\phi\in\mathcal D_K(\Omega).
\]
This shows that $T$ has order $0$.

Now take a sequence of bounded open sets $\{\Omega_i\}_{i\in\mathbb N}$ such that
\[
\Omega_i\Subset\Omega_{i+1},
\qquad
\bigcup_{i=1}^\infty\Omega_i=\Omega ,
\]
and define $K_i:=\overline{\Omega_i}$. Then each $K_i$ is compact and $K_i\subset\Omega$.

Since $T$ has order $0$, for each $i$ there exists $c_i>0$ such that
\begin{equation}\label{eq:cota_local}
|T(\phi)|\le c_i\|\phi\|_{C^0(K_i)},
\qquad \forall\phi\in\mathcal D_{K_i}(\Omega).
\end{equation}

Let $u\in C_c(\Omega_i)$. Since $\supp(u)$ is compact and contained in $\Omega_i$, there exists a sequence
\[
u_k\in\mathcal D_{K_i}(\Omega)
\]
such that
\[
u_k\to u
\]
uniformly on $K_i$.

By \eqref{eq:cota_local},
\[
|T(u_k-u_\ell)|
\le c_i\|u_k-u_\ell\|_{C^0(K_i)}\to0
\qquad (k,\ell\to\infty).
\]
Consequently, $\{T(u_k)\}$ is Cauchy in $\mathbb R$, and we define
\[
T_i(u):=\lim_{k\to\infty}T(u_k).
\]

We show that this definition does not depend on the approximating sequence. Suppose that
\[
v_k\in\mathcal D_{K_i}(\Omega), \qquad v_k\to u
\]
uniformly on $K_i$. Then
\[
|T(u_k-v_k)|
\le c_i\|u_k-v_k\|_{C^0(K_i)}
\le c_i\bigl(\|u_k-u\|_{C^0(K_i)}+\|u-v_k\|_{C^0(K_i)}\bigr)\to0,
\]
and it follows that
\[
\lim_{k\to\infty}T(u_k)=\lim_{k\to\infty}T(v_k).
\]
Therefore, $T_i(u)$ is well defined.

Now let $u,v\in C_c(\Omega_i)$ and $\alpha,\beta\in\mathbb R$. Take sequences
\[
u_k\to u, \qquad v_k\to v
\]
uniformly on $K_i$, with
\[
u_k,v_k\in\mathcal D_{K_i}(\Omega).
\]
Then
\[
\alpha u_k+\beta v_k\to\alpha u+\beta v
\]
uniformly on $K_i$, and by linearity of $T$ we obtain
\[
T_i(\alpha u+\beta v)
=\lim_{k\to\infty}T(\alpha u_k+\beta v_k)
=\alpha T_i(u)+\beta T_i(v).
\]
Thus, $T_i$ is linear.

Moreover, if $u\ge0$, we can choose an approximating sequence $u_k\ge0$. Then
\[
T(u_k)\ge0,
\]
and passing to the limit gives
\[
T_i(u)\ge0.
\]
Therefore, $T_i$ is positive.

Passing to the limit in \eqref{eq:cota_local} gives
\[
|T_i(u)|\le c_i\|u\|_{C^0(K_i)},
\qquad \forall u\in C_c(\Omega_i),
\]
and consequently $T_i$ is continuous with respect to the uniform norm.

If $u\in C_c(\Omega_i)$, any approximating sequence valid for $T_i$ is also valid for $T_{i+1}$, so
\[
T_{i+1}(u)=T_i(u).
\]
Thus the functionals $T_i$ are compatible and define a linear functional
\[
L\colon C_c(\Omega)\longrightarrow\mathbb R
\]
given by
\[
L(u)=T_i(u)\qquad (u\in C_c(\Omega_i)).
\]

Now let $K\subset\Omega$ be compact. Choose $i$ such that $K\subset\Omega_i$. If $\supp(u)\subset K$, then
\[
|L(u)|=|T_i(u)|\le c_i\|u\|_{C^0(K_i)}=c_i\|u\|_{C^0(K)}.
\]
This shows that $L$ is locally bounded and therefore continuous on $C_c(\Omega)$.

Applying the Riesz representation theorem (Theorem~\ref{teo: representacion de Riesz en Cc}), we obtain a unique positive Radon measure $\mu$ such that
\[
L(u)=\int_\Omega u\,d\mu,
\qquad \forall u\in C_c(\Omega).
\]
Since $L$ extends $T$, we conclude that
\[
T(\phi)=\int_\Omega\phi\,d\mu,
\qquad \forall\phi\in\mathcal D(\Omega).
\]

\item Now suppose that $T$ has order $0$. Then, for every compact set $K\subset\Omega$ there exists $C_K>0$ such that
\[
|T(\phi)|\le C_K\|\phi\|_{C^0(K)},
\qquad \forall\phi\in\mathcal D_K(\Omega).
\]

Repeating the preceding procedure, we obtain mutually compatible continuous linear functionals
\[
T_i\colon C_c(\Omega_i)\longrightarrow\mathbb R
\]
inducing a linear functional
\[
L\colon C_c(\Omega)\longrightarrow\mathbb R
\]
that extends $T$.

Moreover, the same argument shows that $L$ is locally bounded. Since $\Omega$ is locally compact, Hausdorff, and $\sigma$-compact, part (b) of Theorem~\ref{teo: representacion de Riesz en Cc} directly provides two positive Radon measures $\mu_1,\mu_2$ such that
\[
L(u)=\int_\Omega u\,d\mu_1-\int_\Omega u\,d\mu_2,
\qquad \forall u\in C_c(\Omega).
\]
Each integral is finite because $u$ has compact support. The difference of the measures need not be defined on every Borel subset of $\Omega$; the representation is understood in the local sense of Definition~\ref{def:b4-espacios-medida-radon-con-signo}. Therefore,
\[
T(\phi)=\int_\Omega\phi\,d\mu_1-\int_\Omega\phi\,d\mu_2,
\qquad \forall\phi\in\mathcal D(\Omega).
\]

\end{enumerate}

\end{proof}

\section{Distributional derivatives}
The notion of a weak derivative motivated the concept of a distribution because, as we saw earlier, some functions have no weak derivatives but are \textit{sufficiently well behaved} in the sense that they should admit some kind of \textit{generalized} derivative. The concept of a distribution will allow us to extend the notion of a derivative as follows:
\begin{definition}\label{def:distribuciones-y-fourier-derivada-conjunto-abierto-multiindice-esima}\index{distributional derivative}
 Let $\Omega\subseteq \mathbb{R}^{n}$ be an open set and let $T\in \mathcal{D}'(\Omega)$. Given a multi-index $\alpha \in \mathbb{N}_{0}^{n}\setminus\{0\}$, we define the $\alpha$th \textit{distributional} derivative of $T$ by \[D^{\alpha}T(\phi):=(-1)^{|\alpha|}T(D^{\alpha}\phi),\qquad \phi \in\mathcal{D}(\Omega).\]
 For $j\in \mathbb{N}$, the symbol $\nabla^{j}T$ represents the collection of all $\alpha$th distributional derivatives of $T$ with $|\alpha|=j$. \end{definition}
\begin{remark}\label{obs: derivada distribucional es distribucion y su orden}
 Naturally, the $\alpha$th distributional derivative of a distribution $T$ must again be a distribution. To prove this, note that by Theorem~\ref{teo: equivalencias de continuidad en D'(Omega)}, if $T\in \mathcal{D}'(\Omega)$ and $K\subseteq \Omega$ is compact, then there exist $j\in \mathbb{N}_{0}$ and a constant $c_{K}>0$ such that
 \[
 |T(\phi)|\leq c_{K}\|\phi\|_{C^{j}(K)},\qquad \forall \phi\in \mathcal{D}_{K}(\Omega).
 \]
 Then
 \[
 |D^{\alpha}T(\phi)|=|T(D^{\alpha}\phi)|\leq c_{K}\|D^{\alpha}\phi\|_{C^{j}(K)}\leq c_{K}\|\phi\|_{C^{j+|\alpha|}(K)},\qquad \forall \phi\in \mathcal{D}_{K}(\Omega).
 \]
 Again by Theorem~\ref{teo: equivalencias de continuidad en D'(Omega)}, we obtain $D^{\alpha}T\in \mathcal{D}'(\Omega)$.

 Moreover, this estimate allows us to control the order of the distributional derivative. If $T$ has finite order $m=\operatorname{ord}(T)$, taking $j=m$ in the preceding inequality gives
 \[
 \operatorname{ord}(D^{\alpha}T)\leq m+|\alpha|.
 \]
 If $T$ does not have finite order, then $\operatorname{ord}(T)=\infty$ and the preceding inequality remains valid in $\mathbb{N}_{0}\cup\{\infty\}$, since
 \[
 \operatorname{ord}(D^{\alpha}T)\leq \infty = \infty+|\alpha|.
 \]

 In particular, if $u\in L^{1}_{\operatorname{loc}}(\Omega)$, $T_{u}$ is the distribution induced by $u$, and $\alpha$ is a multi-index, then the $\alpha$th weak derivative of $u$, $D^{\alpha}u\in L^{1}_{\operatorname{loc}}(\Omega)$ (if it exists), has the property that its induced distribution $T_{D^{\alpha}u}\in \mathcal{D}'(\Omega)$ is the $\alpha$th distributional derivative of $T_{u}$, that is,
 \[
 D^{\alpha}T_{u}=T_{D^{\alpha}u}.
 \]
 If $\phi\in\mathcal{D}(\Omega)$, then by definition of distributional derivative
 \[
 D^{\alpha}T_{u}(\phi)
 =(-1)^{|\alpha|}T_{u}(D^{\alpha}\phi)
 =(-1)^{|\alpha|}\int_{\Omega}u\,D^{\alpha}\phi.
 \]
 Since $D^{\alpha}u$ is the weak derivative of $u$, we have
 \[
 \int_{\Omega}u\,D^{\alpha}\phi
 =
 (-1)^{|\alpha|}\int_{\Omega}D^{\alpha}u\,\phi,
 \]
 and therefore
 \[
 D^{\alpha}T_{u}(\phi)
 =
 \int_{\Omega}D^{\alpha}u\,\phi
 =
 T_{D^{\alpha}u}(\phi).
 \]
 Since this holds for every $\phi\in\mathcal{D}(\Omega)$, we conclude that
 \[
 D^{\alpha}T_{u}=T_{D^{\alpha}u}.
 \]
\end{remark}
 Directional derivatives extend similarly:
 \begin{definition}\label{def:distribuciones-y-fourier-derivada-conjunto-abierto-vector-esima}\index{higher-order distributional directional derivative}
 Let $\Omega\subseteq \mathbb{R}^{n}$ be an open set and let $T\in \mathcal{D}'(\Omega)$. Given a vector $\nu \in \mathbb{R}^{n}\setminus\{0\}$, we define the $k$th directional derivative of $T$ in the direction of $\nu$ by \[\frac{\partial^{k}T}{\partial \nu^{k}}(\phi):=(-1)^{k}T\left(\frac{\partial^{k}\phi}{\partial \nu ^{k}}\right),\qquad \phi\in \mathcal{D}(\Omega).\] In particular, if $u\in L^{1}_{\operatorname{loc}}(\Omega)$, $\nu \in \mathbb{R}^{n}\setminus \{0\}$, and $k\in \mathbb{N}$, then the $k$th weak or distributional directional derivative of $u$ in the direction $\nu$ is the distribution \[\displaystyle\frac{\partial^{k}T_{u}}{\partial \nu^{k}}\in \mathcal{D}'(\Omega).\]
 \end{definition}
 Other differential operators, such as the Laplacian, admit extensions to the distributional sense.
\begin{example}\label{ej:distribuciones-y-fourier-conjunto-abierto-laplaciano-sentido-distribucional-define}
Let $\Omega\subseteq \mathbb{R}^{n}$ be an open set and let $u\in L^{1}_{\operatorname{loc}}(\Omega)$. The Laplacian of $u$ in the distributional sense is defined as the distribution
\[
T_{\Delta u}:=\displaystyle\sum_{i=1}^{n}\partial_{i}^{2}T_{u}.
\]
In other words, for each $\phi\in\mathcal{D}(\Omega)$ we have
\begin{align*}
T_{\Delta u}(\phi)
&=\displaystyle\sum_{i=1}^{n}\partial_{i}^{2}T_{u}(\phi)=\displaystyle\sum_{i=1}^{n}T_{u}(\partial_{i}^{2}\phi)\\
&=\displaystyle\sum_{i=1}^{n}\int_{\Omega}u\,\partial_{i}^{2}\phi=\int_{\Omega}u\,\Delta\phi .
\end{align*}

We say that $u$ is \textit{subharmonic} on $\Omega$ if its distributional Laplacian is positive, that is, if
\[
T_{\Delta u}(\phi)\geq 0,
\qquad
\forall \phi\in\mathcal{D}(\Omega)\text{ with }\phi\geq 0.
\]

In this case, $T_{\Delta u}$ is a positive distribution. By Theorem~\ref{teo: distribuciones orden cero y medidas}, every positive distribution is a positive Radon measure. Therefore, there exists a unique positive Radon measure $\mu$ on $\Omega$ such that
\[
T_{\Delta u}(\phi)=\int_{\Omega}\phi\,d\mu,
\qquad
\forall \phi\in\mathcal{D}(\Omega).
\]

Consequently, the distributional Laplacian of a subharmonic function can be identified with a positive Radon measure.
\end{example}
\begin{example}\label{ej:distribuciones-y-fourier-conjunto-abierto-divergencia-sentido-distribucional-define}
Let $\Omega\subseteq \mathbb{R}^{n}$ be an open set and let $F=(F_{1},\dots,F_{n})\in L^{1}_{\operatorname{loc}}(\Omega,\mathbb{R}^{n})$. The divergence of $F$ in the distributional sense is defined as the distribution
\[
T_{\operatorname{div}F}:=\displaystyle\sum_{i=1}^{n}\partial_{i}T_{F_{i}}.
\]
In other words, for each $\phi\in\mathcal{D}(\Omega)$ we have
\begin{align*}
T_{\operatorname{div}F}(\phi)
&=\displaystyle\sum_{i=1}^{n}\partial_{i}T_{F_{i}}(\phi)=-\displaystyle\sum_{i=1}^{n}T_{F_{i}}(\partial_{i}\phi)\\
&=-\displaystyle\sum_{i=1}^{n}\int_{\Omega}F_{i}\,\partial_{i}\phi=-\int_{\Omega}F\cdot\nabla\phi .
\end{align*}

\end{example}
\begin{exercise}\label{ejer:distribuciones-y-fourier-definida-calcule-derivada-sentido-distribucional-encuentre}
 Let $u\colon \mathbb{R}\longrightarrow \mathbb{R}$ be defined by
 \[
 u(x):=\begin{cases}
 \cos(x),&\text{if $-\pi\leq x\leq 0$,}\\
 1-\displaystyle\frac{x}{\pi},&\text{if $0<x\leq \pi$,}\\
 0,&\text{otherwise.}
 \end{cases}
 \]
 \begin{enumerate}
 \item Compute the derivative of $u$ in the distributional sense and find its order.

 \item Let $v:=u\restriction_{(-\pi,\pi)}$. Compute the first and second derivatives of $v$ in the distributional sense on $(-\pi,\pi)$ and find their orders.

 \item Let $w\colon \mathbb{R}\longrightarrow \mathbb{R}$ be differentiable on $\mathbb{R}\setminus \{a\}$ and suppose that
 \[
 \lim_{x\to a^{-}}w(x),\qquad \lim_{x\to a^{+}}w(x).
 \]
 exist. Compute the derivative of $w$ in the distributional sense. What can you say about its order? Under what conditions can one conclude that it has order $0$?
 \end{enumerate}
\end{exercise}

\begin{solution}
 \begin{enumerate}
 \item Let $\phi\in \mathcal{D}(\mathbb{R})$. Since $u\in L^1_{\operatorname{loc}}(\mathbb{R})$, the function $u$ induces a distribution
 \[
 T_u\colon \mathcal{D}(\mathbb{R})\longrightarrow \mathbb{R},
 \qquad
 T_u(\phi):=\int_{\mathbb{R}}u(x)\phi(x)\,dx.
 \]
 By definition of distributional derivative,
 \[
 T_u'(\phi):=-T_u(\phi')=-\int_{\mathbb{R}}u(x)\phi'(x)\,dx.
 \]
 Since $u$ is given by different expressions on $[-\pi,0]$ and on $(0,\pi]$, we have
 \[
 T_u'(\phi)
 =
 -\int_{-\pi}^{0}\cos(x)\phi'(x)\,dx
 -
 \int_{0}^{\pi}\left(1-\frac{x}{\pi}\right)\phi'(x)\,dx.
 \]

 We integrate each term by parts. For the first,
 \begin{align*}
 -\int_{-\pi}^{0}\cos(x)\phi'(x)\,dx
 &=
 -\Bigl[\cos(x)\phi(x)\Bigr]_{-\pi}^{0}
 +
 \int_{-\pi}^{0}(-\sin x)\phi(x)\,dx\\
 &=
 -\phi(0)-\phi(-\pi)+\int_{-\pi}^{0}(-\sin x)\phi(x)\,dx,
 \end{align*}
 since $\cos(0)=1$ and $\cos(-\pi)=-1$.

 For the second term,
 \[
-\int_{0}^{\pi}\left(1-\frac{x}{\pi}\right)\phi'(x)\,dx
 =
 -\Bigl[\left(1-\frac{x}{\pi}\right)\phi(x)\Bigr]_{0}^{\pi}
 -\frac{1}{\pi}\int_{0}^{\pi}\phi(x)\,dx
\]
\[
=
 \phi(0)-\frac{1}{\pi}\int_{0}^{\pi}\phi(x)\,dx.
\]

 Adding the two expressions gives
 \begin{align*}
 T_u'(\phi)
 &=
 -\phi(-\pi)
 +
 \int_{-\pi}^{0}(-\sin x)\phi(x)\,dx
 -
 \frac{1}{\pi}\int_{0}^{\pi}\phi(x)\,dx.
 \end{align*}
 This can be written as
 \[
 T_u'(\phi)=T_f(\phi)-\delta_{-\pi}(\phi),
 \]
 where $f\colon \mathbb{R}\longrightarrow\mathbb{R}$ is given by
 \[
 f(x):=
 \begin{cases}
 -\sin x,&\text{if $-\pi<x<0$,}\\
 -\displaystyle\frac{1}{\pi},&\text{if $0<x<\pi$,}\\
 0,&\text{otherwise.}
 \end{cases}
 \]
 Therefore,
 \[
 T_u'=T_f-\delta_{-\pi}.
 \]
 In other words,
 \[
 u'=f-\delta_{-\pi}
 \]
 in the distributional sense.

 We now determine its order. Since $f\in L^1_{\operatorname{loc}}(\mathbb{R})$, the distribution $T_f$ has order $0$. Moreover, $\delta_{-\pi}$ also has order $0$. It follows that $u'$, being a sum of distributions of order $0$, is a distribution of order $0$.

 \item Now consider $v:=u\restriction_{(-\pi,\pi)}$. Then
 \[
 v(x)=
 \begin{cases}
 \cos(x),&\text{if $-\pi<x\leq 0$,}\\
 1-\displaystyle\frac{x}{\pi},&\text{if $0<x<\pi$.}
 \end{cases}
 \]
 Moreover,
 \[
 \lim_{x\to 0^-}v(x)=\cos(0)=1
 \qquad\text{and}\qquad
 \lim_{x\to 0^+}v(x)=1-\frac{0}{\pi}=1,
 \]
 so $v$ is continuous at $0$. In particular, $v\in L^1_{\operatorname{loc}}((-\pi,\pi))$, and therefore it induces a distribution
 \[
 T_v\colon \mathcal{D}((-\pi,\pi))\longrightarrow \mathbb{R},
 \qquad
 T_v(\phi):=\int_{-\pi}^{\pi}v(x)\phi(x)\,dx.
 \]

 We compute its first distributional derivative. Let $\phi\in \mathcal{D}((-\pi,\pi))$. By definition,
 \[
 T_v'(\phi):=-T_v(\phi')=-\int_{-\pi}^{\pi}v(x)\phi'(x)\,dx.
 \]
 Splitting the integral according to the definition of $v$ gives
 \[
 T_v'(\phi)
 =
 -\int_{-\pi}^{0}\cos(x)\phi'(x)\,dx
 -
 \int_{0}^{\pi}\left(1-\frac{x}{\pi}\right)\phi'(x)\,dx.
 \]

 We integrate by parts on each interval. For the first,
 \begin{align*}
 -\int_{-\pi}^{0}\cos(x)\phi'(x)\,dx
 &=
 -\Bigl[\cos(x)\phi(x)\Bigr]_{-\pi}^{0}
 +
 \int_{-\pi}^{0}(-\sin x)\phi(x)\,dx\\
 &=
 -\phi(0)-\phi(-\pi)+\int_{-\pi}^{0}(-\sin x)\phi(x)\,dx.
 \end{align*}
 Now, since $\phi\in \mathcal{D}((-\pi,\pi))$, we have $\phi(-\pi)=0$. Consequently,
 \[
 -\int_{-\pi}^{0}\cos(x)\phi'(x)\,dx
 =
 -\phi(0)+\int_{-\pi}^{0}(-\sin x)\phi(x)\,dx.
 \]

 For the second term,
 \[
-\int_{0}^{\pi}\left(1-\frac{x}{\pi}\right)\phi'(x)\,dx
 =
 -\Bigl[\left(1-\frac{x}{\pi}\right)\phi(x)\Bigr]_{0}^{\pi}
 -\frac{1}{\pi}\int_{0}^{\pi}\phi(x)\,dx
\]
\[
=
 \phi(0)-\frac{1}{\pi}\int_{0}^{\pi}\phi(x)\,dx,
\]
 since $\left(1-\displaystyle\frac{\pi}{\pi}\right)\phi(\pi)=0$.

 Adding the two expressions, the terms at $\phi(0)$ cancel and we obtain
 \[
 T_v'(\phi)
 =
 \int_{-\pi}^{0}(-\sin x)\phi(x)\,dx
 -
 \frac{1}{\pi}\int_{0}^{\pi}\phi(x)\,dx.
 \]
 This can be written as
 \[
 T_v'(\phi)=\int_{-\pi}^{\pi}g(x)\phi(x)\,dx,
 \]
 where
 \[
 g(x):=
 \begin{cases}
 -\sin x,&\text{if $-\pi<x<0$,}\\
 -\displaystyle\frac{1}{\pi},&\text{if $0<x<\pi$.}
 \end{cases}
 \]
 Therefore,
 \[
 T_v'=T_g.
 \]
 Consequently,
 \[
 v'=
 \begin{cases}
 -\sin x,&\text{if $-\pi<x<0$,}\\
 -\displaystyle\frac{1}{\pi},&\text{if $0<x<\pi$,}
 \end{cases}
 \qquad\text{in the distributional sense.}
 \]

 Since $g\in L^1_{\operatorname{loc}}((-\pi,\pi))$, the distribution $T_g$ has order $0$. Consequently, $v'$ is a distribution of order $0$.

 We now compute the second distributional derivative of $v$. Since $T_v'=T_g$, we have
 \[
 T_v''=(T_v')'=(T_g)'.
 \]
 Let $\phi\in \mathcal{D}((-\pi,\pi))$. Then
 \[
 T_v''(\phi)=T_g'(\phi)=-T_g(\phi')
 =
 -\int_{-\pi}^{0}(-\sin x)\phi'(x)\,dx
 +\frac{1}{\pi}\int_{0}^{\pi}\phi'(x)\,dx.
 \]

 We integrate each term by parts. For the first,
 \begin{align*}
 -\int_{-\pi}^{0}(-\sin x)\phi'(x)\,dx
 &=
 -\Bigl[(-\sin x)\phi(x)\Bigr]_{-\pi}^{0}
 +
 \int_{-\pi}^{0}(-\cos x)\phi(x)\,dx.
 \end{align*}
 Since $\sin(0)=0$ and $\sin(-\pi)=0$, the boundary term vanishes, so
 \[
 -\int_{-\pi}^{0}(-\sin x)\phi'(x)\,dx
 =
 \int_{-\pi}^{0}(-\cos x)\phi(x)\,dx.
 \]

 For the second term,
 \[
 \frac{1}{\pi}\int_{0}^{\pi}\phi'(x)\,dx
 =
 \frac{1}{\pi}\bigl[\phi(x)\bigr]_{0}^{\pi}
 =
 -\frac{1}{\pi}\phi(0),
 \]
 since $\phi(\pi)=0$.

 Substituting into the expression for $T_v''(\phi)$ gives
 \[
 T_v''(\phi)
 =
 \int_{-\pi}^{0}(-\cos x)\phi(x)\,dx
 -
 \frac{1}{\pi}\phi(0).
 \]
 This means that
 \[
 T_v''=T_h-\frac{1}{\pi}\delta_0,
 \]
 where
 \[
 h(x):=
 \begin{cases}
 -\cos x,&\text{if $-\pi<x<0$,}\\
 0,&\text{if $0<x<\pi$.}
 \end{cases}
 \]
 In other words,
 \[
 v''=
 \begin{cases}
 -\cos x,&\text{if $-\pi<x<0$,}\\
 0,&\text{if $0<x<\pi$,}
 \end{cases}
 -\frac{1}{\pi}\delta_0
 \]
 in the distributional sense.

 We now examine its order. Since $h\in L^1_{\operatorname{loc}}((-\pi,\pi))$, the distribution $T_h$ has order $0$. Moreover, $\delta_0$ has order $0$, so $v''$ is a distribution of order $0$.

 \item Now let $w\colon \mathbb{R}\longrightarrow \mathbb{R}$ be a function differentiable on $\mathbb{R}\setminus\{a\}$, and suppose that the one-sided limits
 \[
 \lim_{x\to a^-}w(x)
 \qquad\text{and}\qquad
 \lim_{x\to a^+}w(x).
 \]
 exist. We denote these limits by
 \[
 w(a^-):=\lim_{x\to a^-}w(x),
 \qquad
 w(a^+):=\lim_{x\to a^+}w(x).
 \]
 By hypothesis, both real numbers are well defined.

 Since $w$ is differentiable on $\mathbb{R}\setminus\{a\}$, it is in particular continuous on $(-\infty,a)$ and on $(a,\infty)$. Moreover, since the finite one-sided limits at $a$ exist, the function $w$ is locally bounded near $a$. It follows that
 \[
 w\in L^1_{\operatorname{loc}}(\mathbb{R}),
 \]
 and therefore $w$ induces a distribution
 \[
 T_w\colon \mathcal{D}(\mathbb{R})\longrightarrow \mathbb{R},
 \qquad
 T_w(\phi):=\int_{\mathbb{R}}w(x)\phi(x)\,dx.
 \]

 We compute its distributional derivative. Let $\phi\in\mathcal{D}(\mathbb{R})$. By definition,
 \[
 T_w'(\phi):=-T_w(\phi')=-\int_{\mathbb{R}}w(x)\phi'(x)\,dx.
 \]
 We split the integral over the intervals $(-\infty,a)$ and $(a,\infty)$:
 \[
 T_w'(\phi)
 =
 -\int_{-\infty}^{a}w(x)\phi'(x)\,dx
 -
 \int_{a}^{\infty}w(x)\phi'(x)\,dx.
 \]

 Suppose also that the classical derivative of $w$ on $\mathbb{R}\setminus\{a\}$, denoted by $w'_{\mathrm{cl}}$, belongs to $L^1_{\operatorname{loc}}(\mathbb{R})$ after being extended arbitrarily at the point $a$. We can then integrate by parts on each of the two preceding intervals. We obtain
 \begin{align*}
 -\int_{-\infty}^{a}w(x)\phi'(x)\,dx
 &=
 -\Bigl[w(x)\phi(x)\Bigr]_{-\infty}^{a}
 +
 \int_{-\infty}^{a}w'_{\mathrm{cl}}(x)\phi(x)\,dx\\
 &=
 -w(a^-)\phi(a)
 +
 \int_{-\infty}^{a}w'_{\mathrm{cl}}(x)\phi(x)\,dx,
 \end{align*}
 since $\phi$ has compact support.

 Similarly,
 \begin{align*}
 -\int_{a}^{\infty}w(x)\phi'(x)\,dx
 &=
 -\Bigl[w(x)\phi(x)\Bigr]_{a}^{\infty}
 +
 \int_{a}^{\infty}w'_{\mathrm{cl}}(x)\phi(x)\,dx\\
 &=
 w(a^+)\phi(a)
 +
 \int_{a}^{\infty}w'_{\mathrm{cl}}(x)\phi(x)\,dx.
 \end{align*}

 Adding the two equalities gives
 \[
 T_w'(\phi)
 =
 \int_{-\infty}^{a}w'_{\mathrm{cl}}(x)\phi(x)\,dx
 +
 \int_{a}^{\infty}w'_{\mathrm{cl}}(x)\phi(x)\,dx
 +
 \bigl(w(a^+)-w(a^-)\bigr)\phi(a).
 \]
 Since a point has measure zero, we can write
 \[
 \int_{-\infty}^{a}w'_{\mathrm{cl}}(x)\phi(x)\,dx
 +
 \int_{a}^{\infty}w'_{\mathrm{cl}}(x)\phi(x)\,dx
 =
 \int_{\mathbb{R}}w'_{\mathrm{cl}}(x)\phi(x)\,dx.
 \]
 Consequently,
 \[
 T_w'(\phi)
 =
 \int_{\mathbb{R}}w'_{\mathrm{cl}}(x)\phi(x)\,dx
 +
 \bigl(w(a^+)-w(a^-)\bigr)\phi(a),
 \]
 that is,
 \[
 T_w'=T_{w'_{\mathrm{cl}}}+\bigl(w(a^+)-w(a^-)\bigr)\delta_a.
 \]
 In other words,
 \[
 w'=w'_{\mathrm{cl}}+\bigl(w(a^+)-w(a^-)\bigr)\delta_a
 \]
 in the distributional sense.

 We now consider the order of $T_w'$. Since $T_w$ is a distribution of order $0$, Remark~\ref{obs: derivada distribucional es distribucion y su orden} implies that $T_w'$ is a distribution of order at most $1$.

 If also $w'_{\mathrm{cl}}\in L^1_{\operatorname{loc}}(\mathbb{R})$, then $T_{w'_{\mathrm{cl}}}$ has order $0$, and since $\delta_a$ also has order $0$, we conclude that
 \[
 T_w'=T_{w'_{\mathrm{cl}}}+\bigl(w(a^+)-w(a^-)\bigr)\delta_a
 \]
 is a distribution of order $0$.

 In particular, this occurs if the two branches of $w$, extended at $a$ with the values $w(a^-)$ and $w(a^+)$, are absolutely continuous on $[b,a]$ and $[a,c]$ for any $b<a<c$. The fundamental theorem of calculus for absolutely continuous functions gives integrability of $|w'_{\mathrm{cl}}|$ on $(b,a)$ and on $(a,c)$. Every compact subset of $\mathbb R$ is contained in some interval $[b,c]$ of this type, so $w'_{\mathrm{cl}}\in L^1_{\operatorname{loc}}(\mathbb R)$. This condition controls the derivative as one approaches $a$; absolute continuity on compact intervals separated from $a$ is insufficient for this.
 \end{enumerate}
\end{solution}
After studying several examples of distributions obtained as derivatives, it is natural to ask to what extent every distribution admits such a representation. The following theorem gives a local answer:
\begin{theorem}\label{teo:distribuciones-y-fourier-abierto-compacto-funcion-continua-multiindice-ambos}\index{distribution!local representation as a derivative}
 Let $\Omega\subseteq \mathbb{R}^{n}$ be open and let $T\in \mathcal{D}'(\Omega)$. Then for every compact set $K\subseteq \Omega$ there exist a continuous function $u\colon \Omega\longrightarrow \mathbb{R}$ and a multi-index $\alpha$ (both depending on $K$) such that
 \[
 T(\phi)=(-1)^{|\alpha|}\int_{\Omega}uD^{\alpha}\phi ,\qquad
 \forall \phi\in \mathcal{D}_{K}(\Omega).
 \]
\end{theorem}

\begin{proof}
 Let $K\subseteq\Omega$ be a fixed compact set. Since $K$ is compact and $\Omega$ is open, there exists an open set $U_0$ such that $K\subseteq U_0\Subset\Omega$. By Proposition~\ref{funciones flan}, there exists $\chi\in C^\infty(\Omega)$ such that $\chi=1$ on $K$ and $\supp(\chi)\subseteq U_0$. In particular, $\chi\in\mathcal D(\Omega)$. Define $S\in\mathcal D'(\mathbb R^n)$ by
 \[
 S(\psi):=T\bigl((\chi\psi)\restriction_\Omega\bigr),
 \qquad \psi\in\mathcal D(\mathbb R^n).
 \]
 Multiplication by $\chi$ and restriction to $\Omega$ are continuous on the corresponding test function spaces; therefore, $S$ is a distribution. Moreover, $\supp(S)\subseteq\supp(\chi)$.

 Choose a closed cube
 \[
 Q=[a_1,b_1]\times\cdots\times[a_n,b_n]\subseteq\mathbb R^n
 \]
 such that $\supp(\chi)\subseteq\operatorname{Int}(Q)$. This choice does not require $Q$ to be contained in $\Omega$. If $\phi\in\mathcal D_K(\Omega)$ and $\widetilde\phi$ denotes its extension by zero to $\mathbb R^n$, then $\widetilde\phi\in\mathcal D_Q(\mathbb R^n)$ and
 \begin{equation}\label{eq:localizacion-distribucion-continua}
 S(\widetilde\phi)=T(\chi\phi)=T(\phi).
 \end{equation}

 Denote by $\beta:=(1,\dots,1)$ the multi-index whose components are all equal to $1$.

 Now let $\psi\in \mathcal{D}_{Q}(\mathbb{R}^n)$. For $x\in Q$ define the possibly degenerate rectangle
 \[
 Q(x):=\{y\in Q \mid a_i\le y_i\le x_i,\ i=1,\dots,n\}.
 \]
 Since $\psi\in\mathcal D_Q(\mathbb R^n)$, its support is contained in the closed cube $Q$. Since $\psi$ is of class $C^\infty$ and vanishes identically on $\mathbb{R}^n \setminus Q$, both $\psi$ and all its partial derivatives must vanish on the boundary $\partial Q$. This ensures that when Theorem~\ref{teo:b5-fundamental-calculo-riemann-banach} is applied in each variable, integrating from the lower endpoint $a_i$, the boundary terms vanish.

 Theorem~\ref{teo:b5-fundamental-calculo-riemann-banach}, applied in the first variable, gives
 \[
 \psi(x)
 =
 \psi(a_1,x_2,\dots,x_n)
 {}+\int_{a_1}^{x_1}\partial_1\psi(y_1,x_2,\dots,x_n)\,dy_1
 =
 \int_{a_1}^{x_1}\partial_1\psi(y_1,x_2,\dots,x_n)\,dy_1.
 \]
 Since $\partial_1\psi$ vanishes on $\partial Q$, it is in particular zero when the second coordinate is exactly $a_2$. This allows us to apply Theorem~\ref{teo:b5-fundamental-calculo-riemann-banach} to the integrand with respect to the second variable without boundary terms appearing. Repeating this argument successively in each of the $n$ variables---and using at each step the fact that the intermediate derivatives also vanish on the boundary of $Q$---gives the integral identity
 \begin{equation}\label{eq:identidad_integral_beta}
 \psi(x)=\int_{Q(x)}D^\beta\psi(y)\,dy,
 \qquad \forall x\in Q.
 \end{equation}

 Since $S\in\mathcal D'(\mathbb R^n)$, the restriction of $S$ to $\mathcal D_Q(\mathbb R^n)$ is continuous by Theorem~\ref{teo: equivalencias de continuidad en D'(Omega)}. Therefore, there exist an integer $j\ge0$ and a constant $C_Q'>0$ such that
 \[
 |S(\psi)|\le C_Q'\|\psi\|_{C^{j}(Q)},
 \qquad \forall\psi\in\mathcal D_Q(\mathbb R^n),
 \]
 where
 \[
 \|\psi\|_{C^{j}(Q)}
 =
 \max_{|\alpha|\le j}\|D^\alpha\psi\|_{C^0(Q)}.
 \]

 Now let $\alpha$ be a multi-index with $|\alpha|\le j$ and set $\gamma:=(j+1)\beta-\alpha$. Since $\alpha_i\le j$, we have $\gamma_i\ge1$ for every $i$. Applying the fundamental theorem of calculus $\gamma_i$ times in the variable $i$ and using the fact that all intermediate derivatives vanish when $x_i=a_i$ gives
 \[
 D^\alpha\psi(x)
 =
 \int_{a_1}^{x_1}\!\cdots\!\int_{a_n}^{x_n}
 \left(\prod_{i=1}^n
 \frac{(x_i-y_i)^{\gamma_i-1}}{(\gamma_i-1)!}\right)
 D^{(j+1)\beta}\psi(y)\,dy_n\cdots dy_1.
 \]
 Indeed, for a single variable this is the repeated integration formula, proved by induction on $\gamma_i$; the displayed formula results from applying it successively in the $n$ variables and using Fubini. Since $0\le x_i-y_i\le b_i-a_i$ on the integration region, taking absolute values and enlarging that region to all of $Q$ gives
 \[
 |D^\alpha\psi(x)|
 \le
 C_{Q,\alpha} \int_Q |D^{(j+1)\beta}\psi(y)|\,dy,
 \]
 where the constant
 \[
 C_{Q,\alpha}
 :=
 \prod_{i=1}^n
 \frac{(b_i-a_i)^{\gamma_i-1}}{(\gamma_i-1)!}
 \]
 depends only on the side lengths of the cube $Q$ and the multi-index $\gamma = (j+1)\beta - \alpha$. Taking the supremum over $x\in Q$, the maximum over $|\alpha|\le j$, and defining
 \[
 C_Q := \max_{|\alpha|\le j} C_{Q,\alpha},
 \]
 we deduce that
 \[
 \|\psi\|_{C^{j}(Q)}
 \le
 C_Q \int_Q |D^{(j+1)\beta}\psi(y)|\,dy.
 \]
 Consequently, absorbing $C_Q$ into a new global constant $C:=C_Q'C_Q$, we obtain the bound
 \begin{equation}\label{eq:cota_distribucion}
 |S(\psi)|
 \le
 C
 \int_Q |D^{(j+1)\beta}\psi(y)|\,dy,
 \qquad \forall\psi\in\mathcal D_Q(\mathbb R^n).
 \end{equation}

 Now define the linear operator
 \[
 L\colon \mathcal D_Q(\mathbb R^n)\longrightarrow L^1(Q),
 \qquad
 L(\psi)=D^{(j+1)\beta}\psi.
 \]
 The operator $L$ is injective. Indeed, suppose that $D^{(j+1)\beta}\psi=0$. We will prove by descending induction that $D^{q\beta}\psi=0$ for $q=j,j-1,\dots,0$. If $q=j$, identity \eqref{eq:identidad_integral_beta}, applied to $D^{j\beta}\psi$, gives $D^{j\beta}\psi=0$ on $Q$. If $0\leq q<j$ and we already know that $D^{(q+1)\beta}\psi=0$, the same identity, now applied to $D^{q\beta}\psi$, gives $D^{q\beta}\psi=0$. Upon reaching $q=0$ we obtain $\psi=0$ on $Q$; the argument also directly covers the case $j=0$. All applications are legitimate because the derivatives of $\psi$ belong to $\mathcal D_Q(\mathbb R^n)$ and hence vanish on $\partial Q$. Outside $Q$ the function $\psi$ is already zero, so $\psi=0$ on $\mathbb R^n$.

 Let $Y:=L(\mathcal D_Q(\mathbb R^n))$ be the image of $L$. Define a linear functional $T_1\colon Y\longrightarrow\mathbb R$ by $T_1(L(\psi)):=S(\psi)$. Injectivity of $L$ ensures that $T_1$ is well defined. Moreover, by \eqref{eq:cota_distribucion},
 \[
 |T_1(w)|
 \le
 C\|w\|_{L^1(Q)},
 \qquad \forall w\in Y.
 \]

 By the analytic form of the Hahn--Banach theorem (Theorem~\ref{teo: hahn--banach, forma analitica}), we can extend $T_1$ to a continuous linear functional on all of $L^1(Q)$. Since the dual space $(L^1(Q))^*$ is isometrically identified with $L^\infty(Q)$ with respect to Lebesgue measure, there exists a function $v\in L^\infty(Q)$ such that
 \[
 T_1(\psi)=\int_Q v\psi,
 \qquad \forall \psi\in L^1(Q).
 \]
 In particular, for $\psi\in\mathcal D_Q(\mathbb R^n)$ we obtain
 \[
 S(\psi)=\int_Q vD^{(j+1)\beta}\psi.
 \]

 Extend $v$ by zero outside $Q$ and define
 \[
 u_0(x)
 =
 \int_{-\infty}^{x_1}\cdots
 \int_{-\infty}^{x_n}
 v(y_1,\dots,y_n)\,dy_n\cdots dy_1.
 \]
 Since $v$ is bounded and has compact support, the function $u_0$ is well defined. Moreover, if $x,x'\in\mathbb R^n$, the symmetric difference between the regions $\{y\in Q:y\leq x\}$ and $\{y\in Q:y\leq x'\}$ is contained in the union of $n$ strips whose measures tend to zero as $x'\to x$. Consequently,
 \[
 |u_0(x')-u_0(x)|
 \leq
 \|v\|_{L^\infty(Q)}
 \lambda_n\bigl(\{y\in Q:y\leq x'\}\mathbin{\triangle}
                 \{y\in Q:y\leq x\}\bigr)
 \longrightarrow0,
 \]
 and $u_0$ is continuous on $\mathbb R^n$. We check that it satisfies the distributional equation $D^\beta u_0=v$. For any test function $\varphi \in \mathcal{D}(\mathbb{R}^n)$, we evaluate the dual action:
 \begin{align*}
 \langle D^\beta u_0, \varphi \rangle
 &=
 (-1)^{|\beta|} \int_{\mathbb{R}^n} u_0(x)\, D^\beta\varphi(x)\,dx \\
 &=
 (-1)^n \int_{\mathbb{R}^n} \left( \int_{\{y \le x\}} v(y_1,\dots,y_n)\,dy_1\cdots dy_n \right) D^\beta\varphi(x)\,dx,
 \end{align*}
 where the notation $y \le x$ means $y_i \le x_i$ for every $i$. By Theorem~\ref{Fubini}, we interchange the order of integration; for fixed $y$, the variable $x$ ranges over the domain $\{x \ge y\}$:
 \[
 \langle D^\beta u_0, \varphi \rangle
 =
 (-1)^n \int_{\mathbb{R}^n} v(y_1,\dots,y_n) \left( \int_{\{x \ge y\}} D^\beta\varphi(x)\,dx \right) dy_1\cdots dy_n.
 \]
 Integrating $D^\beta\varphi = \partial_1\cdots\partial_n\varphi$ in each variable $x_i$ from $y_i$ to $\infty$ gives exactly $(-1)^n \varphi(y_1,\dots,y_n)$ because $\varphi$ has compact support. Therefore,
 \[
 \langle D^\beta u_0, \varphi \rangle
 =
 (-1)^n \int_{\mathbb{R}^n} v(y_1,\dots,y_n)\,(-1)^n\,\varphi(y_1,\dots,y_n)\,dy_1\cdots dy_n
 =
 \langle v, \varphi \rangle.
 \]
 This proves that $D^\beta u_0 = v$ in the sense of distributions.

 Substituting this fact into the preceding expression for $S(\psi)$ and using the definition of distributional derivative, we obtain, for every $\psi\in\mathcal D_Q(\mathbb R^n)$,
 \[
 S(\psi)
 =
 \langle D^\beta u_0,D^{(j+1)\beta}\psi\rangle
 =
 (-1)^{|\beta|}
 \int_{\mathbb R^n}u_0\,D^{(j+2)\beta}\psi
 =
 (-1)^n\int_{\mathbb R^n}u_0\,D^{(j+2)\beta}\psi.
 \]

 Now define the multi-index $\alpha:=(j+2)\beta$, so that $|\alpha| = (j+2)n$. To adjust the sign, define the continuous function
 \[
 \widehat u:=(-1)^{n-|\alpha|}u_0=(-1)^{(j+1)n}u_0,
 \]
 where in the second equality we used $n - |\alpha| = -(j+1)n$ and $(-1)^{-(j+1)n} = (-1)^{(j+1)n}$ (since an integer and its negative have the same parity). Let $u:=\widehat u\restriction_\Omega$. With this choice, we have $(-1)^{|\alpha|}\widehat u=(-1)^nu_0$. Finally, if $\phi\in\mathcal D_K(\Omega)$, then, by \eqref{eq:localizacion-distribucion-continua} and because $\supp(D^\alpha\widetilde\phi)\subseteq K\subseteq\Omega$,
 \[
 T(\phi)
 =
 S(\widetilde\phi)
 =
 (-1)^{|\alpha|}
 \int_{\mathbb R^n}\widehat u\,D^\alpha\widetilde\phi
 =
 (-1)^{|\alpha|}
 \int_{\Omega}u\,D^\alpha\phi.
 \]
 This completes the proof.
\end{proof}

\section{Rapidly decreasing functions and the Schwartz class}
\label{sec:clase-schwartz}

Compactly supported test functions are suitable for studying local phenomena, but compact support is not preserved by the Fourier transform. To analyze a function and its frequencies simultaneously, it is helpful to replace this condition with one controlling behavior at infinity. The Schwartz class consists of smooth functions whose derivatives decrease faster than any power. Its topology records all these controls at once and will be the natural test function space for tempered distributions. The development in this section follows \cite[Chapter~10, \S\S10.4--10.8]{Leoni2017}, with the unitary Fourier normalization to be fixed below.

We work with complex-valued functions. Real-valued functions are viewed through their canonical inclusion into the complex-valued ones; whenever a real result is needed, it will suffice to restrict the resulting operators.

\begin{definition}[Schwartz class]\label{def:clase-schwartz}\index{Schwartz class}\index{rapidly decreasing function}
The \textit{Schwartz class} on $\mathbb R^n$ is
\[
\mathcal S(\mathbb R^n)
:=
\left\{
\varphi\in C^\infty(\mathbb R^n,\mathbb C)
\;\middle|\;
p_{\alpha,\beta}(\varphi)<\infty
\text{ for any }\alpha,\beta\in\mathbb N_0^n
\right\},
\]
where $p_{\alpha,\beta}(\varphi):=\displaystyle\sup_{x\in\mathbb R^n}|x^\alpha D^\beta\varphi(x)|$. For $N\in\mathbb N_0$ we also write
\[
q_N(\varphi)
:=
\max_{|\alpha|+|\beta|\leq N}p_{\alpha,\beta}(\varphi).
\]
We equip $\mathcal S(\mathbb R^n)$ with the locally convex topology generated by the countable family $(q_N)_{N\in\mathbb N_0}$.
\end{definition}

The preceding condition is equivalent to saying that, for each $m\in\mathbb N_0$ and each multi-index $\beta$, the function $(1+\|x\|)^mD^\beta\varphi(x)$ is bounded. In particular, $D^\beta\varphi(x)\to0$ faster than $\|x\|^{-m}$ for any $m$. The function $x\mapsto e^{-\|x\|^2}$ belongs to $\mathcal S(\mathbb R^n)$ and does not have compact support, whereas $\mathcal D(\mathbb R^n)\subseteq\mathcal S(\mathbb R^n)$.

\begin{proposition}[Equivalent seminorms]\label{prop:seminormas-equivalentes-schwartz}
For $N\in\mathbb N_0$ define
\[
r_N(\varphi)
:=
\max_{|\beta|\leq N}
\sup_{x\in\mathbb R^n}(1+\|x\|)^N|D^\beta\varphi(x)|.
\]
The families $(q_N)$ and $(r_N)$ generate the same topology on $\mathcal S(\mathbb R^n)$. More precisely, for each $N$ there exist $M\in\mathbb N_0$ and positive constants $c_N,C_N$ such that $q_N\leq c_Nr_N$ and $r_N\leq C_Nq_M$.
\end{proposition}

\begin{proof}
If $|\alpha|+|\beta|\leq N$, then $|x^\alpha D^\beta\varphi(x)|
\leq(1+\|x\|)^N|D^\beta\varphi(x)|$, giving $q_N\leq r_N$. For the reverse inequality, we use the fact that there exists a constant $C_N>0$ such that
\[
(1+\|x\|)^N
\leq
C_N\displaystyle\sum_{|\alpha|\leq N}|x^\alpha|,
\qquad x\in\mathbb R^n.
\]
This estimate is obtained by expanding $(1+|x_1|+\cdots+|x_n|)^N$ and observing that $\|x\|\leq|x_1|+\cdots+|x_n|$. Therefore, $r_N(\varphi)\leq C_Nq_{2N}(\varphi)$.
\end{proof}

\begin{theorem}[Fréchet structure of the Schwartz class]\label{teo:schwartz-frechet}\index{Frechet space@Fréchet space!Schwartz class}
The space $\mathcal S(\mathbb R^n)$ is a Fréchet space. A sequence $(\varphi_j)$ converges to $\varphi$ in $\mathcal S(\mathbb R^n)$ if and only if $q_N(\varphi_j-\varphi)\to0$ for every $N\in\mathbb N_0$.
\end{theorem}

\begin{proof}
By Proposition~\ref{prop: metrizacion por seminormas}, the countable family $(q_N)$ induces an invariant metric compatible with the topology. It remains to check completeness. Let $(\varphi_j)$ be a Cauchy sequence. For each pair of multi-indices $\alpha,\beta$, the sequence $(x^\alpha D^\beta\varphi_j)_j$ is Cauchy in the Banach space $C_b(\mathbb R^n)$ and converges uniformly to a continuous function $g_{\alpha,\beta}$.

Set $\varphi:=g_{0,0}$. For $i\in\{1,\dots,n\}$, $x\in\mathbb R^n$, and $t\in\mathbb R$, Theorem~\ref{teo:b5-fundamental-calculo-riemann-banach} gives
\[
\varphi_j(x+te_i)-\varphi_j(x)
=
\int_0^tD_i\varphi_j(x+se_i)\,ds.
\]
Uniform convergence allows us to pass to the limit inside the integral and obtain
\[
\varphi(x+te_i)-\varphi(x)
=
\int_0^tg_{0,e_i}(x+se_i)\,ds.
\]
Thus, $D_i\varphi=g_{0,e_i}$. Repeating the argument gives
$\varphi\in C^\infty(\mathbb R^n)$ and $D^\beta\varphi=g_{0,\beta}$
for every $\beta$. At each fixed point $x\in\mathbb R^n$,
\[
x^\alpha D^\beta\varphi(x)
=x^\alpha\lim_{j\to\infty}D^\beta\varphi_j(x)
=\lim_{j\to\infty}x^\alpha D^\beta\varphi_j(x)
=g_{\alpha,\beta}(x).
\]
Since $g_{\alpha,\beta}$ is bounded, this identity proves
$\varphi\in\mathcal S(\mathbb R^n)$. The uniform convergence of the
weighted sequences, already established, then gives
\[
p_{\alpha,\beta}(\varphi_j-\varphi)
=\|x^\alpha D^\beta\varphi_j-g_{\alpha,\beta}\|_{C_b(\mathbb R^n)}
\longrightarrow0.
\]
For each $N$, the maximum defining $q_N$ contains finitely many
seminorms; hence $q_N(\varphi_j-\varphi)\to0$.

\end{proof}

\begin{proposition}[Continuous operations on $\mathcal S$]\label{prop:operaciones-continuas-schwartz}
Let $\alpha\in\mathbb N_0^n$, $a\in\mathbb R^n$, and $A\in\operatorname{GL}(n,\mathbb R)$. The maps
\[
\varphi\longmapsto D^\alpha\varphi,
\qquad
\varphi\longmapsto x^\alpha\varphi,
\qquad
\varphi\longmapsto\varphi(\,\cdot-a),
\qquad
\varphi\longmapsto\varphi\circ A
\]
are continuous linear endomorphisms of $\mathcal S(\mathbb R^n)$. Multiplication $\mathcal S(\mathbb R^n)\times\mathcal S(\mathbb R^n)
\longrightarrow\mathcal S(\mathbb R^n)$ is bilinear and continuous.
\end{proposition}

\begin{proof}
The first two assertions follow directly from the definition. For translations, we use $\|x\|\leq\|x-a\|+\|a\|$ and the binomial formula. If $A\in\operatorname{GL}(n,\mathbb R)$, then
\[
\|A^{-1}\|^{-1}\|x\|
\leq\|Ax\|
\leq\|A\|\,\|x\|,
\qquad x\in\mathbb R^n.
\]
The chain rule and these inequalities show that, for each pair of multi-indices $\alpha,\beta$, there exists a constant $C_{\alpha,\beta,A}>0$, depending only on $\alpha$, $\beta$, $A$, and the fixed Euclidean norm, such that
\[
p_{\alpha,\beta}(\varphi\circ A)
\leq C_{\alpha,\beta,A}
\max_{\substack{|\delta|=|\alpha|\\|\gamma|=|\beta|}}
p_{\delta,\gamma}(\varphi).
\]
Indeed, $\|x\|^{|\alpha|}\leq
\|A^{-1}\|^{|\alpha|}\|Ax\|^{|\alpha|}$, the power $\|Ax\|^{|\alpha|}$ is bounded by a finite sum of monomials of degree $|\alpha|$, and each derivative of $\varphi\circ A$ is a finite linear combination of derivatives of $\varphi$ of order $|\beta|$, with coefficients depending only on $A$. Finally, the Leibniz rule gives
\[
x^\alpha D^\beta(\varphi\psi)
=
\displaystyle\sum_{\gamma\leq\beta}\binom{\beta}{\gamma}
(x^\alpha D^\gamma\varphi)D^{\beta-\gamma}\psi,
\]
and each summand is controlled by a seminorm of $\varphi$ and a seminorm of $\psi$.
\end{proof}

\begin{proposition}[Density]\label{prop:densidad-D-en-S-y-S-en-Lp}
The space $\mathcal D(\mathbb R^n)$ is dense in $\mathcal S(\mathbb R^n)$. Moreover, for every $1\leq p<\infty$, the inclusions
\[
\mathcal D(\mathbb R^n)\subseteq\mathcal S(\mathbb R^n)\subseteq L^p(\mathbb R^n)
\]
have dense image and the inclusion $\mathcal S(\mathbb R^n)\hookrightarrow L^p(\mathbb R^n)$ is continuous.
\end{proposition}

\begin{proof}
Take $\eta\in C_c^\infty(\mathbb R^n)$ with $0\leq\eta\leq1$
and $\eta=1$ on $\overline B_{\mathrm{euc}}(0,1)$, and set
$\eta_R(x):=\eta(\displaystyle\frac{x}{R})$ for $R\geq1$.
If $\varphi\in\mathcal S(\mathbb R^n)$, then
$\eta_R\varphi\in\mathcal D(\mathbb R^n)$. For $|\beta|\leq N$,
the Leibniz rule gives
\[
 D^\beta((1-\eta_R)\varphi)
 =(1-\eta_R)D^\beta\varphi
 -\sum_{\substack{0<\gamma\leq\beta}}
   \binom{\beta}{\gamma}R^{-|\gamma|}
   (D^\gamma\eta)(\frac{x}{R})D^{\beta-\gamma}\varphi.
\]
All terms vanish when $\|x\|<R$. On their support,
$1+\|x\|\geq1+R$ and $R^{-|\gamma|}\leq1$. The definition of
$r_{N+1}$ therefore implies
\[
 (1+\|x\|)^N
 \bigl|D^\beta((1-\eta_R)\varphi)(x)\bigr|
 \leq\frac{r_{N+1}(\varphi)}{1+R}
 \left(1+\sum_{\substack{0<\gamma\leq\beta}}
       \binom{\beta}{\gamma}\|D^\gamma\eta\|_{L^\infty(\mathbb R^n)}
 \right).
\]
Taking the suprema and the finite maximum over $|\beta|\leq N$
gives a constant $C_{N,\eta}$ independent of $R$ such that
\[
 q_N((1-\eta_R)\varphi)
 \leq r_N((1-\eta_R)\varphi)
 \leq\frac{C_{N,\eta}}{1+R}r_{N+1}(\varphi)
 \longrightarrow0.
\]
The final limit follows because $r_{N+1}(\varphi)$ is finite and
$\displaystyle\frac{1}{1+R}\to0$. This proves density in
$\mathcal S(\mathbb R^n)$.

If $m>\displaystyle\frac{n}{p}$, then
\[
\|\varphi\|_{L^p(\mathbb R^n)}^p
\leq
r_m(\varphi)^p
\int_{\mathbb R^n}(1+\|x\|)^{-mp}\,dx,
\]
so the inclusion into $L^p(\mathbb R^n)$ is continuous. Density follows from $\mathcal D(\mathbb R^n)\subseteq\mathcal S(\mathbb R^n)$ and Proposition~\ref{densidad de Cc en Lp}.
\end{proof}

\section{Tempered distributions}
\label{sec:distribuciones-temperadas}

The topology of $\mathcal S$ makes it possible to impose controlled growth on a distribution. Unlike a general distribution, a tempered distribution can be Fourier transformed because the transform preserves the Schwartz class.

\begin{definition}[Tempered distribution]\label{def:distribucion-temperada}\index{distribution!tempered}\index{space of tempered distributions}
\glsadd{distribucion-temperada}
A \textit{tempered distribution} is a continuous linear functional $T\colon \mathcal S(\mathbb R^n)\longrightarrow\mathbb C$. The continuous dual of $\mathcal S(\mathbb R^n)$ is denoted by $\mathcal S'(\mathbb R^n)$ and is equipped, unless otherwise stated, with the weak-$*$ topology $\sigma(\mathcal S'(\mathbb R^n),\mathcal S(\mathbb R^n))$.
\end{definition}

\begin{proposition}[Continuity criterion]\label{prop:criterio-continuidad-distribucion-temperada}
A linear functional $T\colon \mathcal S(\mathbb R^n)\longrightarrow\mathbb C$ belongs to $\mathcal S'(\mathbb R^n)$ if and only if there exist $N\in\mathbb N_0$ and $C>0$ such that
\[
|\langle T,\varphi\rangle|
\leq
Cq_N(\varphi)
\qquad
\text{for every }\varphi\in\mathcal S(\mathbb R^n).
\]
\end{proposition}

\begin{proof}
Sufficiency is immediate. Suppose that $T$ is continuous at the origin. There exists a basic neighborhood
\[
U
=
\left\{\varphi\in\mathcal S(\mathbb R^n)\mid
q_{N_1}(\varphi)<\varepsilon_1,\dots,q_{N_k}(\varphi)<\varepsilon_k\right\}
\]
such that $|T(\varphi)|<1$ for $\varphi\in U$. Since the seminorms $q_N$ are increasing, if $\displaystyle N:=\displaystyle\max_jN_j$ and $\displaystyle \varepsilon:=\min_j\varepsilon_j$, then $q_N(\varphi)<\varepsilon$ implies $\varphi\in U$. For $q_N(\varphi)>0$, we apply the preceding argument to $\displaystyle\frac{\varepsilon\varphi}{2q_N(\varphi)}$ and obtain $|T(\varphi)|\leq2\varepsilon^{-1}q_N(\varphi)$. If $q_N(\varphi)=0$, then $\varphi=0$.
\end{proof}

\begin{proposition}[Natural inclusions]\label{prop:inclusiones-Lp-Sprima}
For every $1\leq p\leq\infty$ there is a continuous injective linear embedding
\[
L^p(\mathbb R^n)\hookrightarrow\mathcal S'(\mathbb R^n),
\qquad
f\longmapsto T_f,
\qquad
\langle T_f,\varphi\rangle:=\int_{\mathbb R^n}f(x)\varphi(x)\,dx.
\]
Every locally integrable function of at most polynomial growth also defines a tempered distribution. Every distribution with compact support is tempered.
\end{proposition}

\begin{proof}
Let $p'$ be the conjugate exponent. By Hölder's inequality in Proposition~\ref{desigualdad de holder} and continuity of $\mathcal S(\mathbb R^n)\hookrightarrow L^{p'}(\mathbb R^n)$,
\[
|\langle T_f,\varphi\rangle|
\leq
\|f\|_{L^p(\mathbb R^n)}\|\varphi\|_{L^{p'}(\mathbb R^n)}
\leq
C\|f\|_{L^p(\mathbb R^n)}q_N(\varphi)
\]
for some $N$. Injectivity follows from Proposition~\ref{prop unicidad derivada debil prop 14.49}.

If $|f(x)|\leq C(1+\|x\|)^m$ almost everywhere, choose $N>m+n$ and estimate
\[
\int_{\mathbb R^n}|f\varphi|
\leq
C r_N(\varphi)
\int_{\mathbb R^n}(1+\|x\|)^{m-N}\,dx.
\]
For the last assertion, let $T\in\mathcal D'(\mathbb R^n)$ have compact support $K$. Choose $\chi\in\mathcal D(\mathbb R^n)$ equal to one on a neighborhood of $K$ and define, for $\varphi\in\mathcal S(\mathbb R^n)$,
\[
\langle\widetilde T,\varphi\rangle
:=
\langle T,\chi\varphi\rangle.
\]
The definition does not depend on the cutoff: if $\chi_1$ and $\chi_2$ equal one near $K$, then $(\chi_1-\chi_2)\varphi$ vanishes on a neighborhood of $\operatorname{supp}(T)$ and, by definition of the support of a distribution, $\langle T,(\chi_1-\chi_2)\varphi\rangle=0$. By the local finite-order estimate for $T$ on $\operatorname{supp}(\chi)$, there exist $C>0$ and $N\in\mathbb N_0$ such that
\[
|\langle\widetilde T,\varphi\rangle|
\leq
C\max_{|\alpha|\leq N}\|D^\alpha(\chi\varphi)\|_{C^0(\operatorname{supp}(\chi))}
\leq C_\chi q_N(\varphi).
\]
Thus, $\widetilde T\in\mathcal S'(\mathbb R^n)$ and its restriction to $\mathcal D(\mathbb R^n)$ agrees with $T$.
\end{proof}

\begin{example}\label{ej:temperadas-polinomios-medidas}
Every polynomial defines a tempered distribution. More generally, if $\mu$ is a Radon measure for which there exist $C,m>0$ with $|\mu|(B_{\mathrm{euc}}(0,R))\leq C(1+R)^m$ for $R\geq1$, then $\langle T_\mu,\varphi\rangle:=\displaystyle\int\varphi\,d\mu$ defines an element of $\mathcal S'(\mathbb R^n)$. To see this, decompose $\mathbb R^n$ into the unit ball and the annuli $k\leq\|x\|<k+1$; the growth of $\mu$ is absorbed by a seminorm $r_N(\varphi)$ with $N>m+1$.
\end{example}

\begin{definition}[Differentiation and multiplication]\label{def:operaciones-distribuciones-temperadas}
Let $T\in\mathcal S'(\mathbb R^n)$ and $\alpha\in\mathbb N_0^n$. Define
\[
\langle D^\alpha T,\varphi\rangle
:=
(-1)^{|\alpha|}\langle T,D^\alpha\varphi\rangle,
\qquad
\langle x^\alpha T,\varphi\rangle
:=
\langle T,x^\alpha\varphi\rangle.
\]
Both functionals belong to $\mathcal S'(\mathbb R^n)$.
\end{definition}

The final assertion follows from continuity of the operations in Proposition~\ref{prop:operaciones-continuas-schwartz}. The same argument allows us to multiply $T$ by any function $a\in C^\infty(\mathbb R^n)$ whose derivatives have at most polynomial growth: the Leibniz rule shows that $\varphi\mapsto a\varphi$ is continuous on $\mathcal S(\mathbb R^n)$.

\section{Convolution and regularization}
\label{sec:convolucion-temperadas}

Convolution connects regularization by mollifiers with multiplication in frequency. We first recall the inequality ensuring its existence in the spaces $L^p(\mathbb R^n)$.

\begin{theorem}[Young's inequality for convolutions]\label{teo:young-convoluciones}\index{Young's inequality!convolutions}
Let $1\leq p,q,r\leq\infty$ satisfy $\displaystyle\frac{1}{p}+\displaystyle\frac{1}{q}=1+\displaystyle\frac{1}{r}$. If $f\in L^p(\mathbb R^n)$ and $g\in L^q(\mathbb R^n)$, then $f*g$ is defined almost everywhere, belongs to $L^r(\mathbb R^n)$, and
\[
\|f*g\|_{L^r(\mathbb R^n)}
\leq
\|f\|_{L^p(\mathbb R^n)}\|g\|_{L^q(\mathbb R^n)}.
\]
In particular,
$L^1(\mathbb R^n)*L^p(\mathbb R^n)\subseteq L^p(\mathbb R^n)$
continuously.
\end{theorem}

\begin{proof}
First consider the nonnegative measurable function
\[
F(x):=\int_{\mathbb R^n}|f(x-y)g(y)|\,dy,
\]
allowing the value $+\infty$. If $r=\infty$, Hölder gives $F(x)\leq\|f\|_{L^p}\|g\|_{L^q}$ for every $x$. If $p=1$, then $r=q$. Apply Theorem~\ref{teo:minkowski-integral} first to bounded truncations of $|f|$ and $|g|$ with bounded support. Monotone convergence and translation invariance give
\[
\|F\|_{L^q(\mathbb R^n)}
\leq\int_{\mathbb R^n}|f(y)|\,
\|g(\,\cdot-y)\|_{L^q(\mathbb R^n)}\,dy
=\|f\|_{L^1(\mathbb R^n)}\|g\|_{L^q(\mathbb R^n)}.
\]
The case $q=1$ follows by exchanging $f$ and $g$, using the change of variable $y\mapsto x-y$ in the nonnegative integral.

In the remaining cases, $1<p,q<r<\infty$. Put
\[
H(x):=\int_{\mathbb R^n}|f(x-y)|^p|g(y)|^q\,dy.
\]
Tonelli and translation invariance yield
\[
\int_{\mathbb R^n}H(x)\,dx
=\|f\|_{L^p(\mathbb R^n)}^p\|g\|_{L^q(\mathbb R^n)}^q<\infty.
\]
Thus $H$ is finite almost everywhere. Write
\[
|f(x-y)g(y)|
=\bigl(|f(x-y)|^p|g(y)|^q\bigr)^{\frac{1}{r}}
|f(x-y)|^{1-\frac{p}{r}}|g(y)|^{1-\frac{q}{r}}.
\]
The reciprocals of $r$, $\displaystyle\frac{pr}{r-p}$, and $\displaystyle\frac{qr}{r-q}$ sum to one. Hölder applied to these three factors gives, almost everywhere,
\[
F(x)\leq H(x)^{\frac{1}{r}}\|f\|_{L^p}^{1-\frac{p}{r}}\|g\|_{L^q}^{1-\frac{q}{r}}.
\]
Raising to the power $r$ and integrating yields
\[
\|F\|_{L^r}^r\leq
\|f\|_{L^p}^{r-p}\|g\|_{L^q}^{r-q}\int_{\mathbb R^n}H(x)\,dx
=\|f\|_{L^p}^r\|g\|_{L^q}^r.
\]
In every case $F$ is finite almost everywhere. The integral defining $f*g$ is therefore absolutely convergent there and satisfies $|f*g|\leq F$. This proves both the asserted existence and inequality.
\end{proof}

\begin{proposition}[Convolution in the Schwartz class]\label{prop:convolucion-schwartz}
If $\varphi,\psi\in\mathcal S(\mathbb R^n)$, then $\varphi*\psi\in\mathcal S(\mathbb R^n)$ and
\[
D^\alpha(\varphi*\psi)
=
(D^\alpha\varphi)*\psi
=
\varphi*(D^\alpha\psi).
\]
The map $(\varphi,\psi)\mapsto\varphi*\psi$ is bilinear and continuous from $\mathcal S(\mathbb R^n)\times\mathcal S(\mathbb R^n)$ to $\mathcal S(\mathbb R^n)$.
\end{proposition}

\begin{proof}
For $x$ in a compact set, the derivatives of $\varphi$ are bounded
and $\psi$ is integrable. The fundamental theorem of calculus bounds
the difference quotients by a constant times $|\psi(y)|$.
Dominated convergence, repeated for each derivative, gives
\[
D^\beta(\varphi*\psi)(x)
=\int_{\mathbb R^n}D^\beta\varphi(x-y)\psi(y)\,dy.
\]
Integration by parts in $y$, with vanishing boundary terms
by rapid decay, transfers all these derivatives to $\psi$.
For the seminorms, the identity
$x^\alpha=((x-y)+y)^\alpha$ gives the exact formula
\[
x^\alpha D^\beta(\varphi*\psi)(x)
=\sum_{\gamma\leq\alpha}\binom{\alpha}{\gamma}
  \int_{\mathbb R^n}(x-y)^{\alpha-\gamma}
       D^\beta\varphi(x-y)\,y^\gamma\psi(y)\,dy.
\]
If $|\alpha|+|\beta|\leq N$, the first factor is bounded
by $q_N(\varphi)$. By the equivalence of seminorms in
Proposition~\ref{prop:seminormas-equivalentes-schwartz},
\[
\|y^\gamma\psi\|_{L^1}
\leq C_Nq_{N+n+1}(\psi)
            \int_{\mathbb R^n}(1+\|y\|)^{-n-1}\,dy.
\]
The integral is finite. Therefore,
\[
q_N(\varphi*\psi)
\leq C_{N,n}q_N(\varphi)q_{N+n+1}(\psi).
\]
This bound proves that convolution belongs to $\mathcal S$ and
that the bilinear map is continuous in its topology.
\end{proof}

For $\varphi\in\mathcal S(\mathbb R^n)$ we write $\check\varphi(x):=\varphi(-x)$.

\begin{definition}[Convolution of a tempered distribution with a Schwartz function]\label{def:convolucion-temperada-schwartz}
Let $T\in\mathcal S'(\mathbb R^n)$ and $\varphi\in\mathcal S(\mathbb R^n)$. Define
\[
(T*\varphi)(x)
:=
\langle T,\varphi(x-\,\cdot)\rangle,
\qquad x\in\mathbb R^n.
\]
\end{definition}

\begin{proposition}[Tempered regularization]\label{prop:regularizacion-temperada}
The function $T*\varphi$ is smooth and
\[
D^\alpha(T*\varphi)
=
T*(D^\alpha\varphi)
=
(D^\alpha T)*\varphi.
\]
For each $\alpha$ there exist $C_\alpha>0$ and $N_\alpha\in\mathbb N_0$ such that $|D^\alpha(T*\varphi)(x)|\leq C_\alpha(1+\|x\|)^{N_\alpha}$. Moreover, if $\rho\in\mathcal S(\mathbb R^n)$ satisfies $\displaystyle\int_{\mathbb R^n}\rho=1$ and $\rho_\varepsilon(x):=\varepsilon^{-n}\rho(\displaystyle\frac{x}{\varepsilon})$, then $T*\rho_\varepsilon\to T$ in $\mathcal S'(\mathbb R^n)$ as $\varepsilon\to0^+$.
\end{proposition}

\begin{proof}
For $x$ in a compact set and $|h|\leq1$, Taylor's formula applied
to each derivative of the kernel gives, in every Schwartz seminorm,
\[
q_N\left(
\frac{\varphi(x+he_i-\,\cdot)-\varphi(x-\,\cdot)}{h}
-D_i\varphi(x-\,\cdot)\right)
\leq C_{N,\varphi,x}|h|,
\]
with a constant uniform on that compact set. Use derivatives of
order up to $N+2$ and
$1+\|y\|\leq(1+\|x-y\|)(1+\|x\|)$ to control the weights.
Continuity of $T$ allows passage to the limit. Repeat the argument
with each derivative to prove smoothness and
\[
D^\alpha(T*\varphi)(x)=\langle T,D^\alpha\varphi(x-\,\cdot)\rangle
=((D^\alpha T)*\varphi)(x).
\]
The second equality uses
$D_y^\alpha\varphi(x-y)=(-1)^{|\alpha|}D^\alpha\varphi(x-y)$.
By Proposition~\ref{prop:criterio-continuidad-distribucion-temperada},
$|\langle T,\zeta\rangle|\leq Cq_N(\zeta)$ for some $N$.
The preceding weight inequality gives
\[
q_N(\varphi(x-\,\cdot))\leq C_{N,\varphi}(1+\|x\|)^N,
\]
and proves polynomial growth.

We justify the pairing identity for arbitrary
$\varphi,\psi\in\mathcal S(\mathbb R^n)$. The function
$V(x)(y):=\varphi(x-y)\psi(x)$ is continuous in $x$ with values in
$\mathcal S$, and for each $N$ it satisfies
\[
q_N(V(x))\leq C_{N,\varphi}(1+\|x\|)^N|\psi(x)|.
\]
The right-hand side is integrable. On a compact cube the Riemann
sums of $V$ are Cauchy in every $q_N$, by uniform continuity;
completeness of $\mathcal S$ gives their limit $H_R$.
Continuity of the evaluation of each derivative identifies
\[
H_R(y)=\int_{[-R,R]^n}\varphi(x-y)\psi(x)\,dx.
\]
The tails satisfy
\[
q_N(H_{R'}-H_R)
\leq C_{N,\varphi}\int_{\mathbb R^n\setminus[-R,R]^n}
(1+\|x\|)^N|\psi(x)|\,dx\longrightarrow0
\]
for $R'\geq R\to\infty$. Thus $H_R$ converges in $\mathcal S$
to $H=\check\varphi*\psi$. Apply $T$ to the Riemann sums and
these limits to obtain
\[
\langle T*\varphi,\psi\rangle
=\int_{\mathbb R^n}\langle T,\varphi(x-\,\cdot)\rangle\psi(x)\,dx
=\langle T,\check\varphi*\psi\rangle.
\]
The first integral is absolutely convergent by the polynomial
growth proved above.

Finally, take $\varphi=\rho_\varepsilon$. Since
$\int_{\mathbb R^n}\check\rho=1$, we have
\[
(\check\rho_\varepsilon*\psi)(x)-\psi(x)
=\int_{\mathbb R^n}\check\rho(z)
\bigl(\psi(x-\varepsilon z)-\psi(x)\bigr)\,dz.
\]
The fundamental theorem of calculus along the segment and
$1+\|x\|\leq(1+\|x-\varepsilon t z\|)(1+\|z\|)$, for
$0\leq t\leq1$ and $0<\varepsilon\leq1$, give
\[
q_N(\check\rho_\varepsilon*\psi-\psi)
\leq C_N\varepsilon q_{N+1}(\psi)
\int_{\mathbb R^n}|\rho(z)|\,\|z\|(1+\|z\|)^N\,dz.
\]
The integral is finite. Hence $\check\rho_\varepsilon*\psi\to\psi$
in $\mathcal S$, and the preceding identity together with continuity
of $T$ proves $T*\rho_\varepsilon\to T$ in $\mathcal S'$.
\end{proof}

\section{The Fourier transform}
\label{sec:transformada-fourier}

The Fourier transform turns derivatives into multiplication and convolutions into products. This correspondence makes it possible to measure regularity through the behavior of a function in frequency. We choose the unitary normalization because it makes the transform an isometry of $L^2(\mathbb R^n)$ without additional constants. For the classical properties, the extension to distributions, and compatibility among the different definitions, we follow the development in \cite[Sections~10.3--10.5]{Brigola2025}, with our normalization in $\mathbb R^n$; the proofs are included to fix all constants and to avoid confusing the integral formula with the extensions by density or duality.

\begin{definition}[Fourier transform on $L^1(\mathbb R^n)$]\label{def:transformada-fourier-L1}\index{Fourier transform}
\glsadd{transformada-fourier}
For $f\in L^1(\mathbb R^n)$ define
\[
\widehat f(\xi)
:=
\mathcal Ff(\xi)
:=
(2\pi)^{-\frac{n}{2}}
\int_{\mathbb R^n}e^{-ix\cdot\xi}f(x)\,dx,
\qquad \xi\in\mathbb R^n.
\]
For $g\in L^1(\mathbb R^n)$, the formal inverse transform is defined by
\[
\mathcal F^{-1}g(x)
:=
(2\pi)^{-\frac{n}{2}}
\int_{\mathbb R^n}e^{ix\cdot\xi}g(\xi)\,d\xi.
\]
\end{definition}

\begin{proposition}[Properties on $L^1(\mathbb R^n)$]\label{prop:fourier-L1-continuidad}
The map $\mathcal F\colon L^1(\mathbb R^n)\longrightarrow C_b(\mathbb R^n)$ is linear and
\[
\|\widehat f\|_{L^\infty(\mathbb R^n)}
\leq
(2\pi)^{-\frac{n}{2}}\|f\|_{L^1(\mathbb R^n)}.
\]
Moreover, $\widehat f$ is uniformly continuous.
\end{proposition}

\begin{proof}
The estimate follows by taking the absolute value inside the integral. If $h\to0$, then
\[
|\widehat f(\xi+h)-\widehat f(\xi)|
\leq
(2\pi)^{-\frac{n}{2}}
\int_{\mathbb R^n}|e^{-ix\cdot h}-1|\,|f(x)|\,dx.
\]
The right-hand side does not depend on $\xi$, converges to zero pointwise in the variable $x$, and is dominated by $2(2\pi)^{-\displaystyle\frac{n}{2}}|f|\in L^1(\mathbb R^n)$. Theorem~\ref{convergencia dominada} proves uniform continuity.
\end{proof}

\begin{proposition}[Fourier calculus on $\mathcal S$]\label{prop:calculo-fourier-schwartz}
If $\varphi\in\mathcal S(\mathbb R^n)$ and $\alpha\in\mathbb N_0^n$, then
\[
\mathcal F(D^\alpha\varphi)(\xi)
=
(i\xi)^\alpha\widehat\varphi(\xi),
\qquad
\mathcal F(x^\alpha\varphi)(\xi)
=
i^{|\alpha|}D^\alpha\widehat\varphi(\xi).
\]
In particular, $\mathcal F$ and $\mathcal F^{-1}$ are continuous linear operators from $\mathcal S(\mathbb R^n)$ to itself.
\end{proposition}

\begin{proof}
The first identity follows by integration by parts; the boundary terms disappear by rapid decrease. The second follows by differentiation under the integral sign. Combining both identities with
\[
\|h\|_{L^1(\mathbb R^n)}
\leq
C_N\sup_{x\in\mathbb R^n}(1+\|x\|)^N|h(x)|,
\qquad N>n,
\]
the Leibniz rule applied to $D^\alpha(x^\beta\varphi)$ and the preceding bound with the weight $(1+\|x\|)^{n+1}$ give explicitly
\[
p_{\alpha,\beta}(\widehat\varphi)
\leq
C_{\alpha,\beta}
q_{n+1+|\alpha|+|\beta|}(\varphi).
\]
Indeed, each summand contains a derivative of order at most $|\alpha|$ and a monomial whose degree, after incorporating the weight, does not exceed $|\beta|+n+1$; the factor $(1+\|x\|)^{-n-1}$ is integrable on $\mathbb R^n$. This proves that $\widehat\varphi\in\mathcal S(\mathbb R^n)$ and establishes continuity. The same calculation, with the sign of the phase reversed, proves the assertion for $\mathcal F^{-1}$.
\end{proof}

\begin{lemma}[Transform of the Gaussian]\label{lem:fourier-gaussiana}
For $a>0$, the function $g_a(x):=e^{-\displaystyle\frac{a\|x\|^2}{2}}$ satisfies
\[
\widehat g_a(\xi)
=
a^{-\frac{n}{2}}e^{-\frac{\|\xi\|^2}{2a}}.
\]
In particular, $g_1$ is fixed by the Fourier transform.
\end{lemma}

\begin{proof}
By Fubini's theorem~\ref{Fubini}, it suffices to treat $n=1$. Let $G(\xi):=\widehat g_a(\xi)$. Proposition~\ref{prop:calculo-fourier-schwartz} and the identity $g_a'=-axg_a$ give $i\xi G=-aiG'$, that is, $G'(\xi)=-\displaystyle\frac{\xi G(\xi)}{a}$. Therefore, $G(\xi)=G(0)e^{-\displaystyle\frac{\xi^2}{2a}}$. To calculate the constant, set $I_a:=\displaystyle\int_\mathbb R e^{-\displaystyle\frac{ax^2}{2}}\,dx$. Tonelli's theorem~\ref{teo:tonelli} and the change to polar coordinates, justified by the change of variables theorem~\ref{teorema de cambio de variable}, give
\[
I_a^2
=
\int_{\mathbb R^2}e^{-\frac{a\|x\|^2}{2}}\,dx
=
2\pi\int_0^\infty e^{-\frac{ar^2}{2}}r\,dr
=
\frac{2\pi}{a}.
\]
Thus, $G(0)=(2\pi)^{-\displaystyle\frac{1}{2}}I_a=a^{-\displaystyle\frac{1}{2}}$. The product of the one-dimensional formulas completes the proof.
\end{proof}

\begin{lemma}[Gaussian approximation to the identity]\label{lem:aproximacion-gaussiana-identidad}
For $\varepsilon>0$, let
\[
k_\varepsilon(x):=(2\pi\varepsilon)^{-\frac{n}{2}}e^{-\frac{\|x\|^2}{2\varepsilon}}.
\]
Then $k_\varepsilon\geq0$, $\displaystyle\int_{\mathbb R^n}k_\varepsilon=1$, and, for every $\delta>0$,
\[
\int_{\{\|x\|\geq\delta\}}k_\varepsilon(x)\,dx\longrightarrow0.
\]
If $1\leq p<\infty$ and $f\in L^p(\mathbb R^n)$, then $k_\varepsilon*f\to f$ in $L^p(\mathbb R^n)$; convergence is uniform when $f$ is bounded and uniformly continuous.
\end{lemma}
\begin{proof}
The mass is computed using Lemma~\ref{lem:fourier-gaussiana}. The change of variable $x=\sqrt\varepsilon z$ and integrability of the Gaussian prove concentration. Moreover, Theorem~\ref{teo:minkowski-integral} gives
\[
\|k_\varepsilon*f-f\|_{L^p(\mathbb R^n)}
\leq
\int_{\mathbb R^n}k_\varepsilon(y)\|\tau_yf-f\|_{L^p(\mathbb R^n)}\,dy.
\]
Split the integral over $\{\|y\|<\delta\}$ and its complement; on the first part use Proposition~\ref{prop:continuidad-traslaciones-Lp}, and on the second use the bound $2\|f\|_{L^p(\mathbb R^n)}$ and concentration. The same argument with the uniform norm proves the final assertion.
\end{proof}

\begin{theorem}[Fourier inversion in the Schwartz class]\label{teo:inversion-fourier-schwartz}\index{Fourier inversion theorem}
For every $\varphi\in\mathcal S(\mathbb R^n)$ we have
\[
\mathcal F^{-1}(\mathcal F\varphi)
=
\varphi
=
\mathcal F(\mathcal F^{-1}\varphi).
\]
Consequently, $\mathcal F\colon \mathcal S(\mathbb R^n)\longrightarrow\mathcal S(\mathbb R^n)$ is a topological isomorphism with inverse $\mathcal F^{-1}$.
\end{theorem}

\begin{proof}
Fix $x\in\mathbb R^n$ and $\varepsilon>0$. Fubini's theorem~\ref{Fubini} and Lemma~\ref{lem:fourier-gaussiana} give
\[
(2\pi)^{-\frac{n}{2}}
\int_{\mathbb R^n}e^{ix\cdot\xi}e^{-\frac{\varepsilon\|\xi\|^2}{2}}\widehat\varphi(\xi)\,d\xi
=
\int_{\mathbb R^n}\varphi(y)k_\varepsilon(x-y)\,dy,
\]
where $k_\varepsilon$ is the function in Lemma~\ref{lem:aproximacion-gaussiana-identidad}. That lemma shows that the right-hand side converges to $\varphi(x)$. On the left-hand side, $e^{-\displaystyle\frac{\varepsilon\|\xi\|^2}{2}}\to1$ and the integrand is dominated by $|\widehat\varphi|\in L^1(\mathbb R^n)$. Theorem~\ref{convergencia dominada} gives $\mathcal F^{-1}\widehat\varphi(x)=\varphi(x)$. The second identity follows by reversing the sign of the phase.
\end{proof}

\begin{corollary}[Inversion in $L^1(\mathbb R^n)$]\label{cor:inversion-fourier-L1}
If $f\in L^1(\mathbb R^n)$ and $\widehat f\in L^1(\mathbb R^n)$, then $f$ agrees almost everywhere with a continuous function vanishing at infinity and
\[
f(x)
=
(2\pi)^{-\frac{n}{2}}
\int_{\mathbb R^n}e^{ix\cdot\xi}\widehat f(\xi)\,d\xi
\]
for almost every $x\in\mathbb R^n$. In particular, the Fourier transform is injective on $L^1(\mathbb R^n)$.
\end{corollary}

\begin{proof}
Let $k_\varepsilon$ be the Gaussian in the preceding proof. By Fubini's theorem~\ref{Fubini},
\[
(f*k_\varepsilon)(x)
=
(2\pi)^{-\frac{n}{2}}
\int_{\mathbb R^n}e^{ix\cdot\xi}e^{-\frac{\varepsilon\|\xi\|^2}{2}}\widehat f(\xi)\,d\xi.
\]
The right-hand side converges uniformly to $\mathcal F^{-1}\widehat f$ by Theorem~\ref{convergencia dominada}. The left-hand side converges to $f$ in $L^1(\mathbb R^n)$ by Lemma~\ref{lem:aproximacion-gaussiana-identidad}. A subsequence converges pointwise almost everywhere, so $f=\mathcal F^{-1}\widehat f$ almost everywhere. Set $h:=\mathcal F^{-1}\widehat f$. By Proposition~\ref{prop:fourier-L1-continuidad}, applied with the opposite sign in the phase, $h$ is uniformly continuous. Moreover, $h=f$ almost everywhere and $f\in L^1(\mathbb R^n)$, so $h\in L^1(\mathbb R^n)$.

These two properties imply that $h$ vanishes at infinity. Indeed, otherwise there would exist $\varepsilon>0$ and a sequence $(x_j)$ with $\|x_j\|\to\infty$ and $|h(x_j)|\geq\varepsilon$. Uniform continuity provides a radius $r>0$, independent of $j$, such that $|h(x)|\geq\displaystyle\frac{\varepsilon}{2}$ if $\|x-x_j\|<r$. Passing to a subsequence, we may assume that the balls $B(x_j,r)$ are pairwise disjoint: after finitely many centers have been chosen, take the next one at a distance greater than $2r$ from all of them. Then
\[
\int_{\mathbb R^n}|h(x)|\,dx
\geq
\sum_{j=1}^{\infty}\int_{B(x_j,r)}|h(x)|\,dx
\geq
\sum_{j=1}^{\infty}\frac{\varepsilon}{2}\lambda_n(B(0,r))
=\infty,
\]
which contradicts integrability of $h$.
\end{proof}

\begin{theorem}[Riemann--Lebesgue lemma]\label{teo:riemann-lebesgue-fourier}\index{Riemann--Lebesgue lemma}
The Fourier transform induces a continuous operator
\[
\mathcal F\colon L^1(\mathbb R^n)\longrightarrow C_0(\mathbb R^n).
\]
In particular, $\widehat f(\xi)\to0$ as $\|\xi\|\to\infty$.
\end{theorem}

\begin{proof}
By Proposition~\ref{prop:densidad-D-en-S-y-S-en-Lp}, there exists $(\varphi_j)\subseteq\mathcal S(\mathbb R^n)$ with $\varphi_j\to f$ in $L^1(\mathbb R^n)$. Proposition~\ref{prop:fourier-L1-continuidad} implies that $\widehat\varphi_j\to\widehat f$ uniformly. Each $\widehat\varphi_j$ belongs to $\mathcal S(\mathbb R^n)\subseteq C_0(\mathbb R^n)$, and $C_0(\mathbb R^n)$ is closed in the uniform norm.
\end{proof}

\begin{proposition}[Convolution, product, and symmetries]\label{prop:identidades-fourier-convolucion}
Let $f,g\in\mathcal S(\mathbb R^n)$. Then
\[
\mathcal F(f*g)
=
(2\pi)^{\frac{n}{2}}\widehat f\,\widehat g,
\qquad
\mathcal F(fg)
=
(2\pi)^{-\frac{n}{2}}\widehat f*\widehat g.
\]
If $a\in\mathbb R^n$, $\lambda\neq0$, and $A\in\operatorname{GL}(n,\mathbb R)$, then
\[
\mathcal F(f(\,\cdot-a))(\xi)=e^{-ia\cdot\xi}\widehat f(\xi),
\qquad
\mathcal F(e^{ia\cdot(\cdot)}f)(\xi)=\widehat f(\xi-a),
\]
and
\[
\mathcal F(f\circ A)(\xi)
=
|\det A|^{-1}\widehat f(A^{-T}\xi).
\]
\end{proposition}

\begin{proof}
Since $f,g\in\mathcal S(\mathbb R^n)$, the function $(x,y)\mapsto f(y)g(x-y)$ is integrable on $\mathbb R^{2n}$ and Fubini's theorem~\ref{Fubini} applies. With $z=x-y$,
\[
\mathcal F(f*g)(\xi)
=(2\pi)^{-\frac{n}{2}}
\int_{\mathbb{R}^{n}}\int_{\mathbb{R}^{n}} e^{-ix\cdot\xi}f(y)g(x-y)\,dy\,dx\] \[=(2\pi)^{-\frac{n}{2}}
\left(\int_{\mathbb{R}^{n}} e^{-iy\cdot\xi}f(y)\,dy\right)
\left(\int_{\mathbb{R}^{n}} e^{-iz\cdot\xi}g(z)\,dz\right)=(2\pi)^{\frac{n}{2}}\widehat f(\xi)\widehat g(\xi).
\]
Applying this identity to $\mathcal F^{-1}\widehat f$ and $\mathcal F^{-1}\widehat g$, and using Theorem~\ref{teo:inversion-fourier-schwartz}, gives $\mathcal F(fg)=(2\pi)^{-\displaystyle\frac{n}{2}}\widehat f*\widehat g$.

For translation, substitute $y=x-a$; for modulation, combine the phase as $e^{-ix\cdot(\xi-a)}$. Finally, with $y=Ax$,
\[
\mathcal F(f\circ A)(\xi)
=(2\pi)^{-\frac{n}{2}}|\det A|^{-1}
\int e^{-iy\cdot A^{-T}\xi}f(y)\,dy,
\]
which is the last formula. Each interchange of order or change of variable is justified by absolute integrability of the Schwartz functions.
\end{proof}

\subsection{The transform on $L^2(\mathbb R^n)$ and on the $L^p(\mathbb R^n)$ scale}

The initial definition of Fourier on integrable functions is insufficient for energy arguments. Plancherel's theorem extends the transform to $L^2(\mathbb R^n)$ and provides the starting point for studying which estimates survive away from the Hilbert space exponent.

\begin{theorem}[Parseval and Plancherel]\label{teo:plancherel-fourier}\index{Plancherel theorem}\index{Parseval's identity}
For $f,g\in\mathcal S(\mathbb R^n)$,
\[
\int_{\mathbb R^n}f(x)\overline{g(x)}\,dx
=
\int_{\mathbb R^n}\widehat f(\xi)\overline{\widehat g(\xi)}\,d\xi.
\]
In particular, $\|\widehat f\|_{L^2(\mathbb R^n)}=\|f\|_{L^2(\mathbb R^n)}$. The Fourier transform extends uniquely to a unitary operator $\mathcal F\colon L^2(\mathbb R^n)\longrightarrow L^2(\mathbb R^n)$, whose inverse is the extension of $\mathcal F^{-1}$.
\end{theorem}

\begin{proof}
Since $f,h\in\mathcal S(\mathbb R^n)\subseteq L^1(\mathbb R^n)$, the corresponding double integral is absolutely convergent. Fubini's theorem~\ref{Fubini} gives the bilinear identity
\[
\int_{\mathbb R^n}\widehat f(\xi)h(\xi)\,d\xi
=
\int_{\mathbb R^n}f(x)\widehat h(x)\,dx
\]
in which both pairings are bilinear. Take $h=\overline{\widehat g}$. By inversion,
\[
\widehat h(x)
=\mathcal F(\overline{\widehat g})(x)
=\overline{\mathcal F^{-1}(\widehat g)(x)}
=\overline{g(x)}.
\]
Substitution of this equality proves Parseval's identity. With $g=f$ we obtain $\|\widehat f\|_{L^2(\mathbb R^n)}=\|f\|_{L^2(\mathbb R^n)}$.

Now let $f\in L^2(\mathbb R^n)$ and choose $f_m\in\mathcal S(\mathbb R^n)$ with $f_m\to f$ in $L^2(\mathbb R^n)$. The isometry shows that $(\widehat f_m)$ is Cauchy; define $\mathcal Ff:=\displaystyle\lim_{m\to\infty}\widehat f_m$ in $L^2(\mathbb R^n)$. If $(\widetilde f_m)$ is another approximation, then
\[
\|\widehat f_m-\widehat{\widetilde f_m}\|_{L^2(\mathbb R^n)}
=\|f_m-\widetilde f_m\|_{L^2(\mathbb R^n)}\longrightarrow0,
\]
so the definition is independent of the sequence and preserves the norm. The same procedure extends $\mathcal F^{-1}$; the inversion identities on $\mathcal S(\mathbb R^n)$ and passage to the limit show that the two extended operators are inverses. In particular, $\mathcal F$ is unitary and surjective. By approximation, this also proves Parseval's identity for every pair of elements of $L^2(\mathbb R^n)$.
\end{proof}

When $f\in L^2(\mathbb R^n)$, the integral expression defining $\widehat f$ may fail to converge. The transform is then understood as the limit in $L^2(\mathbb R^n)$ constructed in the proof. If $f\in L^1(\mathbb R^n)\cap L^2(\mathbb R^n)$, take truncations and regularizations $f_m\in\mathcal S(\mathbb R^n)$ converging simultaneously in $L^1(\mathbb R^n)$ and in $L^2(\mathbb R^n)$. Their transforms converge uniformly to the integral transform by Proposition~\ref{prop:fourier-L1-continuidad}, and converge in $L^2(\mathbb R^n)$ to the Plancherel transform. A subsequence of the second limit converges almost everywhere; by uniqueness of the pointwise limit, the two transforms agree almost everywhere.

\begin{proposition}[$L^p(\mathbb R^n)$ functions as tempered distributions]
\label{prop:Lp-distribucion-temperada-fourier}
For every $1\leq p\leq\infty$, the map
\[
L^p(\mathbb R^n)\longrightarrow\mathcal S'(\mathbb R^n),
\qquad
f\longmapsto T_f,
\qquad
\langle T_f,\varphi\rangle:=\int f(x)\varphi(x)\,dx,
\]
is linear, continuous, and injective. Therefore, $\widehat f:=\widehat{T_f}$ is always defined as a tempered distribution. For $p=2$ it agrees with the Plancherel transform, and for $f\in L^1(\mathbb R^n)$ it agrees with the integral formula.
\end{proposition}

Compatibility with regular distributions is recorded in \cite[Section~10.5]{Brigola2025}.

\begin{proof}
If $p<\infty$, Hölder's inequality in Proposition~\ref{desigualdad de holder} gives
\[
|\langle T_f,\varphi\rangle|
\leq\|f\|_{L^p(\mathbb R^n)}\|\varphi\|_{L^{p'}(\mathbb R^n)}.
\]
Choosing $N>\displaystyle\frac{n}{p'}$, we have
\[
\|\varphi\|_{L^{p'}(\mathbb R^n)}
\leq
\|(1+\|x\|)^{-N}\|_{L^{p'}(\mathbb R^n)}
\sup_{x\in\mathbb R^n}(1+\|x\|)^N|\varphi(x)|.
\]
For $p=\infty$, use the same estimate with $\|\varphi\|_{L^1(\mathbb R^n)}$ and any $N>n$. This proves continuity in the topology of $\mathcal S(\mathbb R^n)$ and therefore $T_f\in\mathcal S'(\mathbb R^n)$. If $T_f=0$, then the integral of $f$ against every function in $C_c^\infty(\mathbb R^n)$ is zero; Proposition~\ref{prop unicidad derivada debil prop 14.49} implies that $f=0$ almost everywhere, proving injectivity.

Compatibility with Plancherel is first checked on $\mathcal S(\mathbb R^n)$ and extended to $L^2(\mathbb R^n)$ by density and continuity. For the integral definition on $L^1(\mathbb R^n)$, the identity follows directly by pairing with a Schwartz function and applying Fubini's theorem~\ref{Fubini}.
\end{proof}

\subsection{Fourier on tempered distributions}

Weak solutions and symbols of polynomial growth require us to go beyond the framework of functions. Since the transform preserves the Schwartz space, duality provides a canonical extension to $\mathcal S'(\mathbb R^n)$.

\begin{definition}[Transform of a tempered distribution]\label{def:fourier-distribucion-temperada}
For $T\in\mathcal S'(\mathbb R^n)$ define $\widehat T\in\mathcal S'(\mathbb R^n)$ and $\mathcal F^{-1}T\in\mathcal S'(\mathbb R^n)$ by
\[
\langle\widehat T,\varphi\rangle
:=
\langle T,\widehat\varphi\rangle,
\qquad
\langle\mathcal F^{-1}T,\varphi\rangle
:=
\langle T,\mathcal F^{-1}\varphi\rangle.
\]
\end{definition}

The definition uses the topological transpose of the transform on $\mathcal S(\mathbb R^n)$. No conjugation is introduced because distributions are viewed as linear functionals.

\begin{theorem}[Fourier on $\mathcal S'$]\label{teo:fourier-isomorfismo-temperadas}
The Fourier transform is a linear, weak-$*$ continuous isomorphism from $\mathcal S'(\mathbb R^n)$ to itself, with inverse $\mathcal F^{-1}$. For every $T\in\mathcal S'(\mathbb R^n)$ and every multi-index $\alpha$,
\[
\mathcal F(D^\alpha T)
=
(i\xi)^\alpha\widehat T,
\qquad
\mathcal F(x^\alpha T)
=
i^{|\alpha|}D^\alpha\widehat T.
\]
If $\varphi\in\mathcal S(\mathbb R^n)$, then
\[
\mathcal F(T*\varphi)
=
(2\pi)^{\frac{n}{2}}\widehat T\,\widehat\varphi.
\]
\end{theorem}

\begin{proof}
The Fourier transform and its inverse are continuous automorphisms
of $\mathcal S$. Their transposes are inverses on $\mathcal S'$.
Moreover, each seminorm of the weak-star topology has the form
$T\mapsto|\langle T,\psi\rangle|$; composing it with Fourier gives
$T\mapsto|\langle T,\widehat\psi\rangle|$. This proves continuity
of the transform and its inverse in that topology.

The identities on $\mathcal S$ give
\[
\begin{aligned}
\langle\widehat{D^\alpha T},\psi\rangle
&=(-1)^{|\alpha|}\langle T,D^\alpha\widehat\psi\rangle
=\langle\widehat T,(i\,\cdot)^\alpha\psi\rangle,\\
\langle\widehat{x^\alpha T},\psi\rangle
&=\langle T,x^\alpha\widehat\psi\rangle
=\langle i^{|\alpha|}D^\alpha\widehat T,\psi\rangle.
\end{aligned}
\]
For the second equality use
$\widehat{D^\alpha\psi}=(i\,\cdot)^\alpha\widehat\psi$ and the
definition of a distributional derivative.

For convolution, Proposition~\ref{prop:regularizacion-temperada}
allows us to write
\[
\langle\widehat{T*\varphi},\psi\rangle
=\langle T*\varphi,\widehat\psi\rangle
=\langle T,\check\varphi*\widehat\psi\rangle.
\]
The Fourier formula for the product of two Schwartz functions,
together with $\widehat{\widehat\varphi}=\check\varphi$, implies
\[
\check\varphi*\widehat\psi
=(2\pi)^{\frac{n}{2}}\widehat{\widehat\varphi\,\psi}.
\]
Therefore,
\[
\langle\widehat{T*\varphi},\psi\rangle
=(2\pi)^{\frac{n}{2}}\langle\widehat T,\widehat\varphi\,\psi\rangle,
\]
which is the convolution identity on $\mathcal S'$.
\end{proof}

\section{Littlewood--Paley dyadic decomposition}
\label{sec:littlewood-paley}

The Fourier transform makes it possible to separate a distribution into frequency bands. Instead of studying the entire spectrum simultaneously, we first isolate the low frequencies and then annuli whose characteristic scale is $2^j$, with $j\in\mathbb N$. This inhomogeneous decomposition preserves the information in the entire distribution and will allow us to measure regularity by observing how each block decreases as frequency increases. It will provide the common foundation for Bessel, Besov, and Triebel--Lizorkin spaces.

\begin{proposition}[Dyadic resolution of unity]\label{prop:resolucion-diadica-unidad}
There exists a radial function $\chi\in C_c^\infty(\mathbb R^n)$ such that $0\leq\chi\leq1$, $\chi=1$ on $\overline B_{\mathrm{euc}}(0,1)$, and $\operatorname{supp}(\chi)\subseteq B_{\mathrm{euc}}(0,2)$. If
\[
\phi_0(\xi):=\chi(\xi),
\qquad
\phi_j(\xi):=\chi(2^{-j}\xi)-\chi(2^{-j+1}\xi),
\quad j\in\mathbb N,
\]
then $\phi_j\in C_c^\infty(\mathbb R^n)$, $\phi_j$ is supported in an annulus $c2^j\leq\|\xi\|\leq C2^j$ for $j\in\mathbb N$, and
\[
\displaystyle\sum_{j\in\mathbb N_0}\phi_j(\xi)=1,
\qquad \xi\in\mathbb R^n.
\]
At each point, only a uniformly bounded number of enlarged annuli intersect.
\end{proposition}

\begin{proof}
The existence of $\chi$ follows from a radial cutoff function. The partial sum is telescoping:
\[
\displaystyle\sum_{j=0}^N\phi_j(\xi)
=
\chi(2^{-N}\xi),
\]
and converges to one. The assertions about supports follow from where $\chi(2^{-j}\xi)$ and $\chi(2^{-j+1}\xi)$ can differ.
\end{proof}

\begin{definition}[Dyadic blocks]\label{def:bloques-diadicos}
Fixing a resolution as above, for $T\in\mathcal S'(\mathbb R^n)$ define
\[
\Delta_jT
:=
\mathcal F^{-1}(\phi_j\widehat T),
\qquad j\in\mathbb N_0.
\]
The distributions $\Delta_jT$ are called the \textit{dyadic blocks} of $T$.
\end{definition}

\begin{theorem}[Littlewood--Paley decomposition]\label{teo:descomposicion-littlewood-paley}\index{Littlewood--Paley decomposition}
For every $T\in\mathcal S'(\mathbb R^n)$,
\[
T
=
\displaystyle\sum_{j\in\mathbb N_0}\Delta_jT
\]
with convergence in $\mathcal S'(\mathbb R^n)$. If $T\in\mathcal S(\mathbb R^n)$, convergence takes place in $\mathcal S(\mathbb R^n)$.
\end{theorem}

\begin{proof}
The partial sum is $\mathcal F^{-1}(\chi(2^{-N}\,\cdot)\widehat T)$. For $A,B\in\mathbb N_0$, use the equivalent Schwartz seminorms
\[
p_{A,B}(\varphi):=\max_{|\beta|\leq B}\sup_{\xi\in\mathbb R^n}
(1+\|\xi\|)^A|D^\beta\varphi(\xi)|.
\]
Set $\chi_N(\xi):=\chi(2^{-N}\xi)$. The Leibniz rule gives
\[
D^\beta((\chi_N-1)\varphi)
=(\chi_N-1)D^\beta\varphi+
\sum_{\substack{0<\gamma\leq\beta}}
\binom\beta\gamma 2^{-N|\gamma|}(D^\gamma\chi)(2^{-N}\xi)
D^{\beta-\gamma}\varphi.
\]
All summands vanish when $\|\xi\|<2^N$. On the complement,
$\displaystyle\frac{1}{1+\|\xi\|}\leq2^{-N}$, and the factors arising from $\chi$ are uniformly bounded for $N\geq0$. Thus there is $C_B>0$, independent of $N$ and $\varphi$, such that
\[
p_{A,B}((\chi_N-1)\varphi)
\leq C_B2^{-N}p_{A+1,B}(\varphi).
\]
The right-hand side tends to zero for each $\varphi\in\mathcal S(\mathbb R^n)$, and uniformly when $\varphi$ ranges over a bounded subset of that space.

Continuity of $\widehat T$ provides $A,B$ and $C_T>0$ such that
$|\langle\widehat T,\psi\rangle|\leq C_Tp_{A,B}(\psi)$. Hence
\[
|\langle(\chi_N-1)\widehat T,\varphi\rangle|
\leq C_TC_B2^{-N}p_{A+1,B}(\varphi).
\]
This proves convergence in $\mathcal S'$, even uniformly on bounded sets of test functions. Continuity of $\mathcal F^{-1}$ gives the first assertion. If $T\in\mathcal S$, applying the first estimate to $\widehat T$ and using continuity of $\mathcal F^{-1}$ on $\mathcal S$ proves the second.
\end{proof}

\begin{theorem}[Bernstein inequalities]\label{teo:desigualdades-bernstein}\index{Bernstein inequality}
Let $1\leq p\leq q\leq\infty$, $m\in\mathbb N_0$, and $R>0$.
\begin{enumerate}[label=(\alph*)]
\item If $f\in L^p(\mathbb R^n)$ and $\operatorname{supp}(\widehat f)\subseteq\overline B_{\mathrm{euc}}(0,R)$, then
\[
\|D^\alpha f\|_{L^q(\mathbb R^n)}
\leq
C R^{|\alpha|+n(\frac{1}{p}-\frac{1}{q})}\|f\|_{L^p(\mathbb R^n)},
\qquad |\alpha|\leq m.
\]
\item If $f\in L^p(\mathbb R^n)$ and $\operatorname{supp}(\widehat f)\subseteq\{\xi\in\mathbb R^n\mid cR\leq\|\xi\|\leq CR\}$, then
\[
C_1R^m\|f\|_{L^p(\mathbb R^n)}
\leq
\displaystyle\sum_{|\alpha|=m}\|D^\alpha f\|_{L^p(\mathbb R^n)}
\leq
C_2R^m\|f\|_{L^p(\mathbb R^n)}.
\]
\end{enumerate}
The constants depend only on the indicated parameters and on $c,C$, not on $R$ or $f$.
\end{theorem}

\begin{proof}
First fix $R=1$. For (a), choose a single function
$\eta_1\in C_c^\infty(\mathbb R^n)$ equal to one on a neighborhood
of the closed unit ball. This choice is independent of $f$.
The identity $\eta_1\widehat f=\widehat f$ and the convolution
formula on $\mathcal S'$ give
\[
f=(2\pi)^{-\frac{n}{2}}(\mathcal F^{-1}\eta_1)*f,
\qquad
D^\alpha f=(2\pi)^{-\frac{n}{2}}(D^\alpha\mathcal F^{-1}\eta_1)*f.
\]
The kernel in the second expression is a Schwartz function.
Young's inequality, with $1+\displaystyle\frac{1}{q}=\displaystyle\frac{1}{p}+\displaystyle\frac{1}{r}$, proves the bound in (a),
including infinite exponents. In particular, the stated distributional
derivatives are represented by $L^q$ functions.

For the lower bound in (b), fix
$\eta_2\in C_c^\infty(\mathbb R^n\setminus\{0\})$ equal to one
on a neighborhood of the whole annulus $c\leq\|\xi\|\leq C$.
If $m=0$, the assertion is immediate. If $m\geq1$, set
\[
D_m(\xi):=\sum_{\ell=1}^n\xi_\ell^{2m},\qquad
m_{m e_\ell}(\xi):=i^{-m}\eta_2(\xi)\frac{\xi_\ell^m}{D_m(\xi)}.
\]
The denominator is positive away from the origin, and these symbols
belong to $C_c^\infty(\mathbb R^n\setminus\{0\})$. Moreover,
\[
\eta_2=\sum_{\ell=1}^n m_{m e_\ell}(i\xi_\ell)^m,
\qquad
f=(2\pi)^{-\frac{n}{2}}\sum_{\ell=1}^n
(\mathcal F^{-1}m_{m e_\ell})*D_\ell^m f.
\]
The first part, using a ball of radius $C$, ensures that
$D_\ell^m f\in L^p$. Young's inequality with the $L^1$ kernels
in the last identity gives the lower bound. The upper bound follows
from (a), also with that ball. All constants are uniform in $f$.

For $R>0$, define $F(x):=f(\displaystyle\frac{x}{R})$.
Since $\widehat F(\xi)=R^n\widehat f(R\xi)$ in $\mathcal S'$,
its spectral support is in the unit ball or the fixed annulus,
respectively. With $\displaystyle\frac{1}{\infty}=0$, change of variables gives
\[
\|F\|_p=R^{\frac{n}{p}}\|f\|_p,\qquad
\|D^\alpha F\|_q=R^{\frac{n}{q}-|\alpha|}\|D^\alpha f\|_q.
\]
Applying the radius-one bounds to $F$ produces the factor
$R^{|\alpha|+n(\displaystyle\frac{1}{p}-\displaystyle\frac{1}{q})}$ in (a) and the factor $R^m$ in (b).
\end{proof}

The finite overlap of the dyadic annuli will be the mechanism allowing us later to compare norms constructed from different resolutions.

\chapter{Interpolation spaces}
\label{cap:espacios-interpolacion}

Intermediate spaces often arise between two spaces measuring different levels of size or regularity. Interpolation makes this idea precise and allows continuity of an operator at intermediate levels to be deduced from its bounds at two endpoints. It also yields compactness results under assumptions to be specified. In this book, the theory will serve three purposes: deriving fractional inequalities from integer-order results, distinguishing the real and complex Sobolev scales, and transferring Euclidean definitions to charts and trivializations. We shall use the functional analysis, integration, and Hardy inequality results from the analysis appendix, as well as the Fourier transform from Chapter~\ref{cap:distribuciones-fourier}. The results of this chapter will be applied to the scales in Chapter~\ref{cap:regularidad-intermedia} and, later, to spaces of functions and sections on manifolds.

We shall develop the Lions--Peetre real method and the Calderón complex method. Our general references will be \cite[Chapters~16--17]{Leoni2017} and \cite{BerghLofstrom1976}; applications to function spaces are compared with \cite[\S1.6]{Triebel2} and \cite[\S2.3]{Schneider2021}.

\begin{semblanzaHistorica}{Interpolating instead of proving anew}
Interpolation theory matured with the recognition that many endpoint estimates contain information about all intermediate levels. The Lions--Peetre methods quantify the decomposition of an element between two spaces, while Calderón's complex method uses holomorphic families. An estimate for an operator on two spaces can thus yield estimates on an entire family of intermediate spaces. In the study of regularity, this allows passage from integer to fractional orders.
\end{semblanzaHistorica}

\section{Compatible pairs of Banach spaces}

Interpolation compares two levels of regularity within a common ambient space. Compatibility ensures that the sum and intersection of the endpoint spaces are well defined and suitable for measuring decompositions.
\label{sec:pares-compatibles-interpolacion}

\begin{definition}[Compatible pair]\label{def:par-compatible-banach}\index{compatible pair of Banach spaces}
A pair $(X_0,X_1)$ of normed spaces is called \textit{compatible} if there are continuous linear embeddings of $X_0$ and $X_1$ into a common Hausdorff topological vector space $X$. We shall identify both spaces with their images in $X$.
\end{definition}

The ambient space serves only to give meaning to equality between elements of $X_0$ and $X_1$. Define
\[
X_0\cap X_1
\quad\text{with}\quad
\|x\|_{X_0\cap X_1}:=\max\{\|x\|_{X_0},\|x\|_{X_1}\},
\]
and
\[
X_0+X_1
:=
\{x_0+x_1\mid x_0\in X_0,\ x_1\in X_1\},
\qquad
\|x\|_{X_0+X_1}
:=
\inf_{\substack{x_0\in X_0,\ x_1\in X_1\\x=x_0+x_1}}
\bigl(\|x_0\|_{X_0}+\|x_1\|_{X_1}\bigr).
\]

\begin{proposition}\label{prop:interseccion-suma-banach}
If $(X_0,X_1)$ is a compatible pair of Banach spaces, then $X_0\cap X_1$ and $X_0+X_1$ are Banach spaces. The inclusions $X_0\cap X_1\hookrightarrow X_j\hookrightarrow
X_0+X_1$, $j=0,1$, are continuous.
\end{proposition}

\begin{proof}
A Cauchy sequence in $X_0\cap X_1$ converges in each $X_j$. The two limits agree in $X$ because the embeddings are continuous and $X$ is Hausdorff; their common limit belongs to the intersection.

For the sum, consider $X_0\times X_1$ with the norm $\|(x_0,x_1)\|:=\|x_0\|_{X_0}+\|x_1\|_{X_1}$ and the subspace
\[
N:=\{(x,-x)\mid x\in X_0\cap X_1\}.
\]
The subspace $N$ is closed: if $(x_k,-x_k)\to(a,b)$ in $X_0\times X_1$, then $x_k\to a$ and $x_k\to-b$ in $X$, whence $a=-b\in X_0\cap X_1$. The map $(x_0,x_1)\mapsto
x_0+x_1$ induces an isometry from the quotient $(X_0\times X_1)/N$ onto $X_0+X_1$. Completeness of quotients of Banach spaces by closed subspaces completes the proof.
\end{proof}

Unless otherwise stated, all compatible pairs used from now on will consist of Banach spaces.

\section{The real method: the \texorpdfstring{$K$}{K}-functional}
\label{sec:metodo-K}

To compare \(x\) with the two endpoint spaces, we decompose it as \(x=x_0+x_1\). The same decomposition may be suitable at one scale and inefficient at another. The following functional chooses the smallest possible cost when the norm of the second component is multiplied by \(t\).

\begin{definition}[$K$-functional]\label{def:funcional-K}\index{K-functional}
Let $(X_0,X_1)$ be a compatible pair. For $x\in X_0+X_1$ and $t>0$, define
\[
K(t,x;X_0,X_1)
:=
\inf_{\substack{x_0\in X_0,\ x_1\in X_1\\x=x_0+x_1}}
\bigl(\|x_0\|_{X_0}+t\|x_1\|_{X_1}\bigr).
\]
When the pair is clear, we shall simply write $K(t,x)$.
\end{definition}

For each $t>0$, $K(t,\cdot)$ is a norm on $X_0+X_1$ and satisfies
\[
\min\{1,t\}\,\|x\|_{X_0+X_1}
\leq K(t,x)
\leq\max\{1,t\}\,\|x\|_{X_0+X_1},
\qquad x\in X_0+X_1.
\]
The parameter $t$ weights the component $x_1$: when $t$ is small, its norm contributes less to the cost. If $X_1$ represents greater regularity than $X_0$, the functional quantifies the cost of separating a more regular part from $x$.

\begin{proposition}[Elementary properties of the $K$-functional]\label{prop:propiedades-funcional-K}
Let $x,y\in X_0+X_1$ and $s,t>0$. Then:
\begin{enumerate}[label=(\alph*)]
\item $K(t,x+y)\leq K(t,x)+K(t,y)$ and $K(t,\lambda x)=|\lambda|K(t,x)$;
\item $t\mapsto K(t,x)$ is increasing and concave;
\item $t\mapsto t^{-1}K(t,x)$ is decreasing;
\item $K(t,x;X_0,X_1)=tK(t^{-1},x;X_1,X_0)$;
\item $\displaystyle K(t,x)\leq\min\{\|x\|_{X_0},t\|x\|_{X_1}\}$ when $x\in X_0\cap X_1$;
\item $\displaystyle \min\{1,t\}\|x\|_{X_0+X_1}\leq K(t,x)\leq\displaystyle\max\{1,t\}\|x\|_{X_0+X_1}$.
\end{enumerate}
\end{proposition}

\begin{proof}
The algebraic properties follow by adding decompositions and taking the infimum. Each function $t\mapsto\|x_0\|_{X_0}+t\|x_1\|_{X_1}$ is affine and increasing; the infimum of this family is increasing and concave. Moreover, $t^{-1}K(t,x)$ is the infimum of the decreasing functions $t\mapsto t^{-1}\|x_0\|_{X_0}+\|x_1\|_{X_1}$, taken over decompositions of $x$, and is therefore decreasing. Symmetry follows by interchanging $x_0$ and $x_1$. The last two estimates follow from the decompositions $x=x+0$, $x=0+x$ and from comparing the coefficients $1,t$ with their minimum and maximum.
\end{proof}

Intermediate regularity requires simultaneous control at all scales \(t\in(0,\infty)\). The factor \(t^{-\theta}\) fixes the intermediate level; the exponent \(q\) determines how the costs are combined. The measure \(\displaystyle\frac{dt}{t}\) assigns the same weight to intervals related by dilation.

\begin{definition}[Real interpolation space]\label{def:espacio-interpolacion-real-K}\index{real interpolation!K-method}
\glsadd{interpolacion-real}
Let $0<\theta<1$ and $1\leq q\leq\infty$. Define
\[
(X_0,X_1)_{\theta,q}
:=
\left\{x\in X_0+X_1\mid\|x\|_{\theta,q;K}<\infty\right\},
\]
where, if $q<\infty$,
\[
\|x\|_{\theta,q;K}
:=
\left(
\int_0^\infty
\bigl(t^{-\theta}K(t,x)\bigr)^q\,\frac{dt}{t}
\right)^{\frac{1}{q}},
\]
and $\|x\|_{\theta,\infty;K}:=\displaystyle\sup_{t\in(0,\infty)}t^{-\theta}K(t,x)$.
\end{definition}

\begin{theorem}[Basic properties of the $K$-method]\label{teo:propiedades-basicas-metodo-K}
If $(X_0,X_1)$ is a compatible pair of Banach spaces, then $(X_0,X_1)_{\theta,q}$ is a Banach space and
\[
X_0\cap X_1
\hookrightarrow
(X_0,X_1)_{\theta,q}
\hookrightarrow
X_0+X_1
\]
continuously. If $X_1\hookrightarrow X_0$, let $C_E$ be the norm of the embedding, let $\displaystyle T_E:=\displaystyle\max\{1,C_E\}$, and set, for $q<\infty$,
\[
N_{\theta,q}(x):=\|x\|_{X_0}+\left(
\int_0^1(t^{-\theta}K(t,x))^q\,\frac{dt}{t}
\right)^{\frac{1}{q}}
\]
and, for $q=\infty$,
\[
N_{\theta,\infty}(x):=\|x\|_{X_0}
+\sup_{t\in(0,1]}t^{-\theta}K(t,x).
\]
Then, if $q<\infty$,
\[
\frac{\|x\|_{\theta,q;K}}{\max\{1,(\theta q)^{-\frac{1}{q}}\}}
\leq N_{\theta,q}(x)
\leq
\bigl(1+(\theta q)^{\frac{1}{q}}T_E^\theta\bigr)\|x\|_{\theta,q;K},
\]
whereas
\[
\|x\|_{\theta,\infty;K}
\leq N_{\theta,\infty}(x)
\leq(1+T_E^\theta)\|x\|_{\theta,\infty;K}.
\]
\end{theorem}

\begin{proof}
If $x\in X_0\cap X_1$, Proposition~\ref{prop:propiedades-funcional-K} gives $\displaystyle K(t,x)\leq\min\{1,t\}\|x\|_{X_0\cap X_1}$. Since
\[
\int_0^1t^{(1-\theta)q}\,\frac{dt}{t}
+
\int_1^\infty t^{-\theta q}\,\frac{dt}{t}
<\infty,
\]
we obtain the first embedding when $q<\infty$. For $q=\infty$, use $\displaystyle \sup_{t\in(0,\infty)}t^{-\theta}\min\{1,t\}=1$ directly. For the second embedding, Proposition~\ref{prop:propiedades-funcional-K} implies $\|x\|_{X_0+X_1}\leq K(t,x)$ for $1\leq t\leq2$. Thus, if $q<\infty$,
\[
\|x\|_{X_0+X_1}
\leq
(\log2)^{-\frac{1}{q}}2^\theta\|x\|_{\theta,q;K};
\]
if $q=\infty$, evaluation at $t=1$ suffices.

Now let $(x_j)$ be Cauchy in $(X_0,X_1)_{\theta,q}$. The second embedding implies $x_j\to
x$ in $X_0+X_1$ for some $x$. For $q<\infty$, continuity of the norm $K(t,\cdot)$ on $X_0+X_1$ and Fatou's lemma~\ref{lem:fatou} give
\[
\|x_j-x\|_{\theta,q;K}^q
\leq
\liminf_{k\to\infty}\|x_j-x_k\|_{\theta,q;K}^q.
\]
The case $q=\infty$ follows by taking the supremum after passing to the limit for each $t$. Thus $x_j\to x$ in the interpolation space.

If $X_1\hookrightarrow X_0$, let $C_E$ be the norm of the embedding and let $\displaystyle T_E:=\displaystyle\max\{1,C_E\}$. Then $X_0+X_1=X_0$ and, for every decomposition $x=x_0+x_1$,
\[
\|x\|_{X_0}
\leq
\|x_0\|_{X_0}+C_E\|x_1\|_{X_1}.
\]
Consequently, if $t\geq T_E$, we have $K(t,x)\geq\|x\|_{X_0}$; the decomposition $x=x+0$ gives the reverse inequality. Hence
\[
\int_{T_E}^\infty
\bigl(t^{-\theta}K(t,x)\bigr)^q\,\frac{dt}{t}
=
\frac{T_E^{-\theta q}}{\theta q}
\|x\|_{X_0}^q.
\]
Moreover, $K(t,x)\leq\|x\|_{X_0}$ for every $t>0$. If $\displaystyle I_0:=\int_0^1(t^{-\theta}K(t,x))^q\,\displaystyle\frac{dt}{t}$, it follows that
\[
 \|x\|_{\theta,q;K}
 \leq\bigl(I_0+(\theta q)^{-1}\|x\|_{X_0}^q\bigr)^{\frac{1}{q}}
 \leq\max\{1,(\theta q)^{-\frac{1}{q}}\}N_{\theta,q}(x).
\]
The tail identity starting at $T_E$ gives, in the opposite direction, $\|x\|_{X_0}\leq(\theta q)^{\displaystyle\frac{1}{q}}T_E^\theta
\|x\|_{\theta,q;K}$, while $I_0^{\displaystyle\frac{1}{q}}\leq\|x\|_{\theta,q;K}$. This proves both bounds for $q<\infty$. For $q=\infty$, the inequality $K(t,x)\leq\|x\|_{X_0}$ shows that $\|x\|_{\theta,\infty;K}\leq N_{\theta,\infty}(x)$, and the equality $K(T_E,x)=\|x\|_{X_0}$ gives $\|x\|_{X_0}\leq T_E^\theta\|x\|_{\theta,\infty;K}$; this yields the remaining bound.
\end{proof}

\begin{lemma}[Dyadic discretization]\label{lem:discretizacion-funcional-K}
For $0<\theta<1$ and $1\leq q\leq\infty$, there are constants $c_{\theta,q},C_{\theta,q}>0$ such that
\[
c_{\theta,q}\left\|
\bigl(2^{-k\theta}K(2^k,x)\bigr)_{k\in\mathbb Z}
\right\|_{\ell^q(\mathbb Z)}
\leq
\|x\|_{\theta,q;K}
\leq
C_{\theta,q}
\left\|
\bigl(2^{-k\theta}K(2^k,x)\bigr)_{k\in\mathbb Z}
\right\|_{\ell^q(\mathbb Z)}.
\]
The constants depend only on $\theta$ and $q$.
\end{lemma}

\begin{proof}
If $2^k\leq t\leq2^{k+1}$, monotonicity of $K$ and $\displaystyle\frac{K(t,x)}{t}$ gives
\[
2^{-\theta}2^{-k\theta}K(2^k,x)
\leq
t^{-\theta}K(t,x)
\leq
2^{1-\theta}2^{-k\theta}K(2^k,x).
\]
Integrate over each interval $[2^k,2^{k+1}]$ with respect to $\displaystyle\frac{dt}{t}$ and sum. The argument with suprema covers $q=\infty$.
\end{proof}

\begin{corollary}[Monotonicity in the second index]\label{cor:monotonia-q-interpolacion-real}
If $1\leq q_0\leq q_1\leq\infty$, then
\[
(X_0,X_1)_{\theta,q_0}
\hookrightarrow
(X_0,X_1)_{\theta,q_1}.
\]
Moreover, $(X_0,X_1)_{\theta,q}=(X_1,X_0)_{1-\theta,q}$ and, for every $x$,
\[
\|x\|_{\theta,q;K,(X_0,X_1)}
=\|x\|_{1-\theta,q;K,(X_1,X_0)}.
\]
\end{corollary}

\begin{proof}
The first assertion follows from Lemma~\ref{lem:discretizacion-funcional-K} and $\ell^{q_0}(\mathbb Z)\hookrightarrow\ell^{q_1}(\mathbb Z)$. For the second, use $K(t,x;X_0,X_1)=tK(t^{-1},x;X_1,X_0)$ and the change of variable $s=t^{-1}$.
\end{proof}

\begin{proposition}[Multiplicative inequality]\label{prop:desigualdad-multiplicativa-interpolacion}
There is $C_{\theta,q}>0$ such that for every $x\in X_0\cap X_1$,
\[
\|x\|_{(X_0,X_1)_{\theta,q}}
\leq
C_{\theta,q}
\|x\|_{X_0}^{1-\theta}\|x\|_{X_1}^{\theta}.
\]
\end{proposition}

\begin{proof}
If $x=0$, the assertion is immediate. Otherwise, set $a:=\|x\|_{X_0}>0$, $b:=\|x\|_{X_1}>0$, and $t_0:=\displaystyle\frac{a}{b}$. If $q<\infty$, Proposition~\ref{prop:propiedades-funcional-K} gives
\[
\begin{aligned}
\|x\|_{\theta,q;K}^q
&\leq b^q\int_0^{t_0}t^{(1-\theta)q}\,\frac{dt}{t}
       +a^q\int_{t_0}^\infty t^{-\theta q}\,\frac{dt}{t}\\
&=\left(\frac{1}{(1-\theta)q}+\frac{1}{\theta q}\right)
  a^{(1-\theta)q}b^{\theta q}.
\end{aligned}
\]
For $q=\infty$, the function $\displaystyle t^{-\theta}\min\{a,tb\}$ attains its maximum at $t=t_0$, where its value is $a^{1-\theta}b^\theta$.
\end{proof}

\section{Interpolation of operators}

The analytic value of the $K$-method lies in transferring known endpoint bounds to intermediate spaces. The following result makes this transfer precise and retains control of the operator norm.
\label{sec:interpolacion-operadores-K}

\begin{theorem}[Interpolation theorem for the $K$-method]\label{teo:interpolacion-operadores-metodo-K}\index{interpolation of operators!K-method}
Let $(X_0,X_1)$ and $(Y_0,Y_1)$ be compatible pairs of Banach spaces. Let $T\colon
X_0+X_1\longrightarrow Y_0+Y_1$ be linear, and suppose that its restrictions $T\colon
X_j\longrightarrow Y_j$ are continuous, with norms $M_j$, $j=0,1$. Then
\[
T\colon (X_0,X_1)_{\theta,q}\longrightarrow(Y_0,Y_1)_{\theta,q}
\]
is continuous and
\[
\|T\|
\leq
M_0^{1-\theta}M_1^\theta.
\]
\end{theorem}

\begin{proof}
If $x=x_0+x_1$, then
\[
K(t,Tx;Y_0,Y_1)
\leq
M_0\|x_0\|_{X_0}+tM_1\|x_1\|_{X_1}.
\]
Taking the infimum,
\[
K(t,Tx;Y_0,Y_1)
\leq
M_0K\left(t\frac{M_1}{M_0},x;X_0,X_1\right).
\]
The conclusion follows from the change of variable $s=\displaystyle\frac{tM_1}{M_0}$. If one of the norms $M_j$ is zero, approximate it by $M_j+\varepsilon$ and let $\varepsilon\to0^+$.
\end{proof}

\begin{corollary}[Interpolation of embeddings]\label{cor:interpolacion-encajes-metodo-K}
If $X_j\hookrightarrow Y_j$ for $j\in\{0,1\}$, then $(X_0,X_1)_{\theta,q}\hookrightarrow(Y_0,Y_1)_{\theta,q}$. If there are linear operators $R\colon X_j\longrightarrow Y_j$ and $S\colon Y_j\longrightarrow X_j$, continuous at both endpoints and satisfying $SR=I$ on $X_0+X_1$, then $S$ induces a retraction and $R$ a coretraction between the interpolated spaces.
\end{corollary}
\begin{proof}
The inclusion \(I\colon X_0+X_1\to Y_0+Y_1\) restricted to each endpoint is a continuous operator \(I\colon X_j\to Y_j\). The operator interpolation theorem, applied to \(I\), provides a continuous operator
\[
I\colon (X_0,X_1)_{\theta,q}
\longrightarrow
(Y_0,Y_1)_{\theta,q},
\]
which is precisely the embedding in the first assertion.

In the second situation, applying the same theorem separately to \(R\) and \(S\) gives continuous operators
\[
R\colon (X_0,X_1)_{\theta,q}\longrightarrow(Y_0,Y_1)_{\theta,q},
\qquad
S\colon (Y_0,Y_1)_{\theta,q}\longrightarrow(X_0,X_1)_{\theta,q}.
\]
For \(x\in(X_0,X_1)_{\theta,q}\subseteq X_0+X_1\), the assumed identity on the sum space implies
\[
S(Rx)=x.
\]
Thus \(S\circ R=I\) on the interpolated space: \(S\) is a retraction, and \(R\) is its corresponding coretraction.
\end{proof}

\begin{theorem}[Compactness into a common target space]\label{teo:interpolacion-operador-compacto-destino-comun}
Let $Y$ be a Banach space and let $T\colon X_0+X_1\longrightarrow Y$ be linear. If $T\colon
X_0\longrightarrow Y$ is compact and $T\colon X_1\longrightarrow Y$ is continuous, then
\[
T\colon (X_0,X_1)_{\theta,q}\longrightarrow Y
\]
is compact for $0<\theta<1$ and $1\leq q\leq\infty$.
\end{theorem}

\begin{proof}
Let $B$ be a bounded subset of $(X_0,X_1)_{\theta,q}$. By Lemma~\ref{lem:discretizacion-funcional-K}, including when $q=\infty$, there is $C>0$ such that $K(t,x)\leq Ct^\theta$ for every $t>0$ and $x\in B$. For each $x\in B$, choose $x=a_{x,t}+b_{x,t}$ with
\[
\|a_{x,t}\|_{X_0}+t\|b_{x,t}\|_{X_1}
\leq
2Ct^\theta.
\]
Given a tolerance $\varepsilon>0$, choose $t$ sufficiently large that $2C\|T\|_{X_1\to Y}t^{\theta-1}<\varepsilon$. The set $\{a_{x,t}\mid x\in B\}$ is bounded in $X_0$ and, by compactness of $T\colon X_0\longrightarrow Y$, its image can be covered by finitely many balls of radius $\varepsilon$. Since
\[
\|Tb_{x,t}\|_Y
\leq2C\|T\|_{X_1\to Y}t^{\theta-1}<\varepsilon,
\]
$T(B)$ admits a finite cover by balls of radius $2\varepsilon$. Thus $T(B)$ is totally bounded; since $Y$ is complete, its closure is compact. This proves the assertion.
\end{proof}

\section{The \texorpdfstring{$J$}{J}-method}
\label{sec:metodo-J}

The $K$-functional describes an element through a single scale-dependent decomposition. The $J$-method instead uses a continuous superposition of elements belonging to both endpoint spaces simultaneously.

\begin{definition}[$J$-functional]\label{def:funcional-J}\index{J-functional}
For $a\in X_0\cap X_1$ and $t>0$, define
\[
J(t,a;X_0,X_1)
:=
\max\{\|a\|_{X_0},t\|a\|_{X_1}\}.
\]
\end{definition}

A function $u\colon (0,\infty)\longrightarrow X_0\cap X_1$ is admissible for $x\in X_0+X_1$ if it is strongly measurable, the Bochner integral $\displaystyle\int_0^\infty u(t)\,\displaystyle\frac{dt}{t}$ converges in $X_0+X_1$, and
\[
x
=
\int_0^\infty u(t)\,\frac{dt}{t}.
\]

\begin{definition}[Space defined by the $J$-method]\label{def:espacio-interpolacion-J}\index{real interpolation!J-method}
Let $0<\theta<1$ and $1\leq q\leq\infty$. Define $(X_0,X_1)_{\theta,q;J}$ as the set of $x\in X_0+X_1$ admitting a representation as above for which
\[
t\longmapsto t^{-\theta}J(t,u(t))
\]
belongs to $L^q((0,\infty),\displaystyle\frac{dt}{t})$. Its norm is the infimum of the $L^q((0,\infty),\displaystyle\frac{dt}{t})$ norm over all admissible representations.
\end{definition}

\begin{lemma}[From a $J$-representation to a $K$-estimate]\label{lem:J-controla-K}
If $x=\displaystyle\int_0^\infty u(t)\,\displaystyle\frac{dt}{t}$ is an admissible representation with $t^{-\theta}J(t,u(t))\in L^q((0,\infty),\displaystyle\frac{dt}{t})$, then for $r>0$,
\[
K(r,x)
\leq
\int_0^r\|u(t)\|_{X_0}\,\frac{dt}{t}
+
r\int_r^\infty\|u(t)\|_{X_1}\,\frac{dt}{t}.
\]
\end{lemma}

\begin{proof}
Split the integral over $(0,r)$ and $(r,\infty)$. If $g(t):=t^{-\theta}J(t,u(t))\in L^q((0,\infty),\displaystyle\frac{dt}{t})$, then
\[
\int_0^r\|u(t)\|_{X_0}\,\frac{dt}{t}
\leq
\int_0^rt^\theta g(t)\,\frac{dt}{t}<\infty
\]
and
\[
\int_r^\infty\|u(t)\|_{X_1}\,\frac{dt}{t}
\leq
\int_r^\infty t^{\theta-1}g(t)\,\frac{dt}{t}<\infty.
\]
Finiteness follows from Proposition~\ref{desigualdad de holder}, or directly if $q=1,\infty$. By the Bochner integrability criterion in Proposition~\ref{prop:b5-criterio-integrabilidad-bochner}, the first part converges in $X_0$ and the second in $X_1$. This splitting is an admissible decomposition for $K(r,x)$, and the triangle inequality for the Bochner integral gives the stated estimate.
\end{proof}

\begin{theorem}[Equivalence of the $K$- and $J$-methods]\label{teo:equivalencia-metodos-K-J}\index{real interpolation!equivalence of the K- and J-methods}
For every compatible pair of Banach spaces, $0<\theta<1$ and $1\leq q\leq\infty$,
\[
(X_0,X_1)_{\theta,q;K}
=
(X_0,X_1)_{\theta,q;J}
\]
and there are constants $c_{\theta,q},C_{\theta,q}>0$, independent of the pair and the element, such that
\[
c_{\theta,q}\|x\|_{\theta,q;K}
\leq\|x\|_{\theta,q;J}
\leq C_{\theta,q}\|x\|_{\theta,q;K}
\]
for every $x$ in the common space.
\end{theorem}

\begin{proof}
First let $x$ be given by a $J$-representation. Lemma~\ref{lem:J-controla-K}, followed by the two integral Hardy inequalities from Theorem~\ref{teo:desigualdades-hardy-integrales},
\[
\left\|r^{-\theta}\int_0^rh(t)\,\frac{dt}{t}\right\|_{L^q((0,\infty),\frac{dr}{r})}
\leq
\frac{1}{\theta}
\|t^{-\theta}h(t)\|_{L^q((0,\infty),\frac{dt}{t})}
\]
and
\[
\left\|r^{1-\theta}\int_r^\infty h(t)\,\frac{dt}{t}\right\|_{L^q((0,\infty),\frac{dr}{r})}
\leq
\frac{1}{1-\theta}
\|t^{1-\theta}h(t)\|_{L^q((0,\infty),\frac{dt}{t})},
\]
gives $\|x\|_{\theta,q;K}\leq C_{\theta,q}\|x\|_{\theta,q;J}$. In the first inequality, we use $h(t)=\|u(t)\|_{X_0}$, and in the second, $h(t)=\|u(t)\|_{X_1}$.

For the reverse inclusion, use a dyadic construction. Choose decompositions $x=a_k+b_k$ such that
\[
\|a_k\|_{X_0}+2^k\|b_k\|_{X_1}
\leq
2K(2^k,x).
\]
Lemma~\ref{lem:discretizacion-funcional-K} gives $K(2^k,x)\leq C2^{k\theta}\|x\|_{\theta,q;K}$. Set $v_k:=a_{k+1}-a_k=b_k-b_{k+1}$. Then
\[
\|a_k\|_{X_0}\longrightarrow0\quad(k\to-\infty),
\qquad
\|b_k\|_{X_1}\longrightarrow0\quad(k\to\infty).
\]
Truncating the sum gives $\displaystyle\sum_{k=-N}^{N}v_k=a_{N+1}-a_{-N}\to x$ in $X_0+X_1$. Therefore, $x=\displaystyle\sum_{k\in\mathbb Z}v_k$ in $X_0+X_1$. Moreover,
\[
J(2^k,v_k)
\leq
C\bigl(K(2^k,x)+K(2^{k+1},x)\bigr).
\]
Define $u(t):=(\log 2)^{-1}v_k$ for $2^k\leq t<2^{k+1}$. For $k\leq0$, use $\|v_k\|_{X_0}\leq C2^{k\theta}\|x\|_{\theta,q;K}$; for $k>0$, use $\|v_k\|_{X_1}\leq C2^{-k(1-\theta)}\|x\|_{\theta,q;K}$. Both geometric series converge. Thus $u$ is integrable in $X_0+X_1$ and $x=\displaystyle\int_0^\infty u(t)\,\displaystyle\frac{dt}{t}$. Lemma~\ref{lem:discretizacion-funcional-K} and the preceding estimate give $\|x\|_{\theta,q;J}\leq C\|x\|_{\theta,q;K}$.
\end{proof}

\begin{lemma}[Discrete characterization of the $J$-method]
\label{lem:caracterizacion-discreta-metodo-J}
Let $0<\theta<1$ and $1\leq q\leq\infty$. An element $x\in X_0+X_1$ belongs to $(X_0,X_1)_{\theta,q}$ if and only if it admits a representation
\[
x=\displaystyle\sum_{k\in\mathbb Z}u_k
\quad\text{in }X_0+X_1,
\qquad u_k\in X_0\cap X_1,
\]
for which
\[
\bigl(2^{-k\theta}J(2^k,u_k;X_0,X_1)\bigr)_{k\in\mathbb Z}
\in\ell^q(\mathbb Z).
\]
Moreover, the interpolation norm is equivalent to the infimum of the norm of this sequence over all such representations. If $q<\infty$, $X_0\cap X_1$ is dense in $(X_0,X_1)_{\theta,q}$.
\end{lemma}

\begin{proof}
First suppose that a discrete representation is given, and write $c_k:=2^{-k\theta}J(2^k,u_k)$. For $k\leq0$,
\[
\|u_k\|_{X_0}\leq2^{k\theta}c_k,
\]
and, for $k>0$,
\[
\|u_k\|_{X_1}\leq2^{-k(1-\theta)}c_k.
\]
Proposition~\ref{desigualdad de holder} for sequences and summation of the geometric series show that the two tails converge in $X_0$ and $X_1$, respectively. Defining $u(t):=(\log 2)^{-1}u_k$ on $[2^k,2^{k+1})$ yields an admissible representation and, for every $t\in[2^k,2^{k+1})$,
\[
\frac{2^{-\theta}}{\log 2}\,
2^{-k\theta}J(2^k,u_k)
\leq
t^{-\theta}J(t,u(t))
\leq
\frac{2}{\log 2}\,
2^{-k\theta}J(2^k,u_k).
\]
This proves one of the norm estimates.

Conversely, given an admissible representation $u(t)$, set
\[
u_k:=\int_{2^k}^{2^{k+1}}u(t)\,\frac{dt}{t}.
\]
Before taking these integrals, set
$g(t):=t^{-\theta}J(t,u(t))\in L^q((0,\infty),\displaystyle\frac{dt}{t})$.
The norm of $X_0\cap X_1$ is $\max\{\|\cdot\|_{X_0},\|\cdot\|_{X_1}\}$.
For $t\in[2^k,2^{k+1}]$,
\[
\|u(t)\|_{X_0\cap X_1}
\leq\max\{1,2^{-k}\}J(t,u(t))
\leq\max\{1,2^{-k}\}2^{(k+1)\theta}g(t).
\]
This interval has measure $\log2$ for $\displaystyle\frac{dt}{t}$. Hölder's inequality,
also for $q=1$ or $q=\infty$, proves integrability of the norm
there. Strong measurability of the representation is in the
intersection, so the Bochner integral exists in $X_0\cap X_1$.
The continuous inclusions into $X_0$, $X_1$, and $X_0+X_1$
transport that same integral to the other spaces. Thus,
$u_k\in X_0\cap X_1$.

Proposition~\ref{prop:b5-integral-bochner-operadores-lineales} allows these integrals to be summed in $X_0+X_1$ and gives $x=\displaystyle\sum_{k\in\mathbb Z}u_k$. The triangle inequality for the Bochner integral and $2^k\leq t$ imply
\[
J(2^k,u_k)
\leq
\int_{2^k}^{2^{k+1}}J(t,u(t))\,\frac{dt}{t}.
\]
Proposition~\ref{desigualdad de holder} on each dyadic interval, followed by Theorem~\ref{teo:tonelli}, yields
\[
\left\|
\bigl(2^{-k\theta}J(2^k,u_k)\bigr)_{k\in\mathbb Z}
\right\|_{\ell^q(\mathbb Z)}
\leq
C_{\theta,q}
\|t^{-\theta}J(t,u(t))\|_{L^q((0,\infty),\frac{dt}{t})}.
\]
The case $q=\infty$ follows by taking suprema. Taking the infimum over representations and using Theorem~\ref{teo:equivalencia-metodos-K-J} completes the equivalence of norms. If $q<\infty$, the finite sums $\displaystyle\sum_{\substack{k\in\mathbb Z\\|k|\leq N}}u_k\in X_0\cap X_1$ converge to $x$ in the interpolation norm, since tails of a sequence in $\ell^q(\mathbb Z)$ tend to zero. This proves the final assertion.
\end{proof}

\section{Reiteration}
\label{sec:reiteracion-interpolacion-real}

Interpolation can be applied more than once. The following result states that interpolating between two already interpolated levels produces no new scales.

\begin{lemma}[Holmstedt's formula]\label{lem:formula-holmstedt}
Let $(X_0,X_1)$ be a compatible pair of Banach spaces, and let $0<\theta_0<\theta_1<1$ and $1\leq
q_0,q_1\leq\infty$. Set $X_{\theta_j}:=(X_0,X_1)_{\theta_j,q_j}$. There are constants $c,C>0$ such that for $x\in X_0+X_1$ and $r>0$,
\[
c\bigl(A_0(r,x)+r^{\theta_1-\theta_0}A_1(r,x)\bigr)
\leq
K(r^{\theta_1-\theta_0},x;X_{\theta_0},X_{\theta_1})
\leq
C\bigl(A_0(r,x)+r^{\theta_1-\theta_0}A_1(r,x)\bigr),
\]
where
\[
A_0(r,x)
:=
\left\|t^{-\theta_0}K(t,x;X_0,X_1)\right\|_{L^{q_0}((0,r),\frac{dt}{t})},
\]
and
\[
A_1(r,x)
:=
\left\|t^{-\theta_1}K(t,x;X_0,X_1)\right\|_{L^{q_1}((r,\infty),\frac{dt}{t})},
\]
with suprema when some $q_j=\infty$. If $x\notin X_{\theta_0}+X_{\theta_1}$, the functional in the middle is interpreted as $+\infty$. The constants $c$ and $C$ depend only on $\theta_0,\theta_1,q_0,q_1$.
\end{lemma}

\begin{proof}
Set $\delta:=\theta_1-\theta_0$ and choose $m\in\mathbb Z$ so that $2^m\leq r<2^{m+1}$. First suppose that $A_0(r,x)+r^\delta A_1(r,x)<\infty$. For each $k\in\mathbb Z$, choose a decomposition $x=a_k+b_k$ such that
\[
\|a_k\|_{X_0}+2^k\|b_k\|_{X_1}
\leq2K(2^k,x;X_0,X_1)
\]
and set $u_k:=a_{k+1}-a_k=b_k-b_{k+1}$. Comparison of $K$ on dyadic intervals gives $K(2^k,x)\leq C2^{k\theta_0}A_0(r,x)$ for $2^{k+1}\leq r$ and $K(2^k,x)\leq C2^{k\theta_1}A_1(r,x)$ for $2^{k-1}\geq r$. Thus $a_k\to0$ in $X_0$ as $k\to-\infty$, and $b_k\to0$ in $X_1$ as $k\to\infty$. Therefore, $x=\displaystyle\sum_{k\in\mathbb Z}u_k$ in $X_0+X_1$, and
\[
\|u_k\|_{X_0}+2^k\|u_k\|_{X_1}
\leq
C\bigl(K(2^k,x)+K(2^{k+1},x)\bigr).
\]
Define
\[
x_0:=\displaystyle\sum_{\substack{k\in\mathbb Z\\k\leq m}}u_k,
\qquad
x_1:=\displaystyle\sum_{\substack{k\in\mathbb Z\\k>m}}u_k.
\]
The discrete characterization of the $J$-method in Lemma~\ref{lem:caracterizacion-discreta-metodo-J}, Lemma~\ref{lem:discretizacion-funcional-K}, and comparison of $r$ with $2^m$ give
\[
\|x_0\|_{X_{\theta_0}}
\leq CA_0(r,x),
\qquad
\|x_1\|_{X_{\theta_1}}
\leq CA_1(r,x).
\]
For the terms near the splitting point, use $K(2^{m+1},x)\leq4K(2^{m-1},x)$ and the interval $[2^{m-1},2^m]\subseteq(0,r)$; for the other tail, use $[2^{m+1},2^{m+2}]\subseteq(r,\infty)$. Thus the constants also control the terms adjacent to $r$. Each series converges in $X_0+X_1$ and defines a $J$-representation of the indicated space; this proves membership even when $q_j=\infty$, without requiring convergence of partial sums in the norm of that space. Since $x=x_0+x_1$, this decomposition proves
\[
K(r^\delta,x;X_{\theta_0},X_{\theta_1})
\leq C\bigl(A_0(r,x)+r^\delta A_1(r,x)\bigr).
\]
If either tail is infinite, this inequality holds automatically.

For the reverse bound, let $x=y_0+y_1$ with $y_j\in X_{\theta_j}$. Lemma~\ref{lem:discretizacion-funcional-K} implies the pointwise estimate
\[
K(t,y_j;X_0,X_1)
\leq C t^{\theta_j}\|y_j\|_{X_{\theta_j}},
\qquad t>0.
\]
The triangle inequality for $K(t,\cdot)$ gives
\[
A_0(r,x)
\leq
C\bigl(\|y_0\|_{X_{\theta_0}}+r^\delta
\|y_1\|_{X_{\theta_1}}\bigr)
\]
and
\[
r^\delta A_1(r,x)
\leq
C\bigl(\|y_0\|_{X_{\theta_0}}+r^\delta
\|y_1\|_{X_{\theta_1}}\bigr).
\]
Here we integrate only the powers $t^{\delta q_0}$ on $(0,r)$ and $t^{-\delta q_1}$ on $(r,\infty)$; if an exponent is infinite, take suprema. Taking the infimum over decompositions $x=y_0+y_1$ yields the remaining inequality.
\end{proof}

\begin{theorem}[Real reiteration theorem]\label{teo:reiteracion-interpolacion-real}\index{reiteration theorem!real interpolation}
Let $0\leq\theta_0<\theta_1\leq1$, $0<\eta<1$, and $1\leq q,q_0,q_1\leq\infty$. If $X_{\theta_j}=(X_0,X_1)_{\theta_j,q_j}$ when $0<\theta_j<1$, while $X_{\theta_0}=X_0$ if $\theta_0=0$ and $X_{\theta_1}=X_1$ if $\theta_1=1$, then
\[
(X_{\theta_0},X_{\theta_1})_{\eta,q}
=
(X_0,X_1)_{\theta,q},
\qquad
\theta:=(1-\eta)\theta_0+\eta\theta_1,
\]
and there are constants $c,C>0$, depending only on $\theta_0,\theta_1,\eta,q,q_0,q_1$, such that
\[
c\|x\|_{(X_0,X_1)_{\theta,q}}
\leq\|x\|_{(X_{\theta_0},X_{\theta_1})_{\eta,q}}
\leq C\|x\|_{(X_0,X_1)_{\theta,q}}
\]
for every $x$ in the common space.
\end{theorem}

\begin{proof}
First assume $0<\theta_0<\theta_1<1$ and set $\delta:=\theta_1-\theta_0$, so that $\theta=\theta_0+\eta\delta$. Write
\[
K_*(\tau,x):=K(\tau,x;X_{\theta_0},X_{\theta_1}),
\qquad
c_k:=2^{-k\theta}K(2^k,x;X_0,X_1).
\]
The proof of Lemma~\ref{lem:discretizacion-funcional-K}, with geometric base $2^\delta$ in place of $2$, provides constants $d_*,D_*>0$ such that
\[
d_*\bigl\|(2^{-k\eta\delta}K_*(2^{k\delta},x))_k\bigr\|_{\ell^q(\mathbb Z)}
\leq\|x\|_{(X_{\theta_0},X_{\theta_1})_{\eta,q}}
\leq D_*\bigl\|(2^{-k\eta\delta}K_*(2^{k\delta},x))_k\bigr\|_{\ell^q(\mathbb Z)}.
\]
By Lemma~\ref{lem:formula-holmstedt}, the sequence between the norm signs is bounded above and below, term by term, by multiples of
\[
B_{0,k}:=2^{-k\eta\delta}A_0(2^k,x),
\qquad
B_{1,k}:=2^{k(1-\eta)\delta}A_1(2^k,x).
\]

Splitting the integrals defining $A_0$ and $A_1$ into the intervals $[2^j,2^{j+1}]$ and using monotonicity of $K$ and $\displaystyle\frac{K(t,x)}{t}$, as in the proof of the discretization lemma, reduces the upper bound to the two sequence operators
\[
\begin{aligned}
(\mathsf H^-_{a,r}c)_k
&:=\left(\sum_{\substack{j\in\mathbb Z\\j\leq k}}
       (2^{-a(k-j)}|c_j|)^r\right)^{\frac{1}{r}},\\
(\mathsf H^+_{a,r}c)_k
&:=\left(\sum_{\substack{j\in\mathbb Z\\j\geq k}}
       (2^{-a(j-k)}|c_j|)^r\right)^{\frac{1}{r}},
\end{aligned}
\]
with suprema if $r=\infty$, taking $(a,r)=(\eta\delta,q_0)$ and $(a,r)=((1-\eta)\delta,q_1)$, respectively. For any $a>0$ and $1\leq q,r\leq\infty$, there is $C_{a,q,r}>0$ such that
\[
\|\mathsf H^-_{a,r}c\|_{\ell^q(\mathbb Z)}
+\|\mathsf H^+_{a,r}c\|_{\ell^q(\mathbb Z)}
\leq C_{a,q,r}\|c\|_{\ell^q(\mathbb Z)}.
\]
Indeed, if $q\geq r<\infty$, raise to the power $r$ and apply Young's inequality in $\ell^{\displaystyle\frac{q}{r}}(\mathbb Z)$ to the geometric kernel $(2^{-arm}\mathbf1_{\{m\geq0\}})_m$. If $q<r$, use
\[
\left(\sum_{j\in\mathbb Z}d_j^r\right)^{\frac{q}{r}}\leq\sum_{j\in\mathbb Z}d_j^q
\]
and sum first over the outer index; the geometric series $\displaystyle \sum_{m\in\mathbb N_0}2^{-aqm}$ appears. For $r=\infty$, the supremum is bounded by the sum of powers $q$, and for $q=\infty$ the assertion follows directly. The tail operator is treated with the reflected kernel. We deduce
\[
\|x\|_{(X_{\theta_0},X_{\theta_1})_{\eta,q}}
\leq C\|(c_k)_k\|_{\ell^q(\mathbb Z)}
\leq C'\|x\|_{(X_0,X_1)_{\theta,q}}.
\]

For the reverse bound, the interval $[2^k,2^{k+1}]$ in the definition of $A_0(2^{k+1},x)$ and monotonicity of $K$ give a constant $d_0=d_0(\theta_0,q_0)>0$ such that
\[
2^{-(k+1)\eta\delta}A_0(2^{k+1},x)
\geq d_0c_k.
\]
This also holds for $q_0=\infty$, interpreting the norm as a supremum. Holmstedt's lower bound and the two discretizations yield
\[
\|x\|_{(X_0,X_1)_{\theta,q}}
\leq C\|x\|_{(X_{\theta_0},X_{\theta_1})_{\eta,q}}.
\]

Now consider $\theta_0=0<\theta_1<1$. We first prove
\begin{equation}
K(r^{\theta_1},x;X_0,X_{\theta_1})
\asymp K(r,x)+r^{\theta_1}A_1(r,x),
\label{eq:holmstedt-extremo-izquierdo}
\end{equation}
where $\asymp$ means that each side is bounded by a positive constant times the other, uniformly in $x$ and $r$.

If the right-hand side is finite and $x\neq0$, choose $x=a+b$ with
\[
\|a\|_{X_0}+r\|b\|_{X_1}\leq2K(r,x).
\]
For $t<r$, use $K(t,b)\leq t\|b\|_{X_1}$; for $t\geq r$, the triangle inequality and monotonicity give
\[
K(t,b)\leq K(t,x)+\|a\|_{X_0}\leq3K(t,x).
\]
Integrating on both intervals, or taking suprema, yields
\[
\|b\|_{X_{\theta_1}}
\leq C\bigl(r^{-\theta_1}K(r,x)+A_1(r,x)\bigr).
\]
This decomposition proves one bound in \eqref{eq:holmstedt-extremo-izquierdo}. For the other, write $x=y_0+y_1$, where $y_0\in X_0$ and $y_1\in X_{\theta_1}$. The pointwise estimate for $K$ and the integral of $t^{-\theta_1q_1}$ over $(r,\infty)$ give
\[
\begin{aligned}
K(r,x)&\leq\|y_0\|_{X_0}
       +Cr^{\theta_1}\|y_1\|_{X_{\theta_1}},\\
r^{\theta_1}A_1(r,x)&\leq C\|y_0\|_{X_0}
       +r^{\theta_1}\|y_1\|_{X_{\theta_1}}.
\end{aligned}
\]
Take the infimum over $y_0,y_1$.

The case $0<\theta_0<\theta_1=1$ follows by applying the preceding formula to the pair $(X_1,X_0)$ with parameter $1-\theta_0$ and scale $r^{-1}$. Indeed, symmetry of the $K$-functional and the substitution $t\mapsto t^{-1}$ transform that formula into
\begin{equation}
K(r^{1-\theta_0},x;X_{\theta_0},X_1)
\asymp A_0(r,x)+r^{-\theta_0}K(r,x).
\label{eq:holmstedt-extremo-derecho}
\end{equation}
Multiply \eqref{eq:holmstedt-extremo-izquierdo} by $r^{-\eta\theta_1}$, or \eqref{eq:holmstedt-extremo-derecho} by $r^{-\eta(1-\theta_0)}$. In both cases, the term containing $K(r,x)$ is $r^{-\theta}K(r,x)$. The other term is controlled by $\mathsf H^+_{(1-\eta)\theta_1,q_1}$ or $\mathsf H^-_{\eta(1-\theta_0),q_0}$, respectively. The estimates already proved for these operators give the upper bound, and the term $r^{-\theta}K(r,x)$ gives the lower bound. If $\theta_0=0$ and $\theta_1=1$, the assertion is the definition.
\end{proof}
\begin{definition}[$\ell^q(I)$-sums of a family of Banach spaces]
\label{def:sumas-ellq-familia-banach}
Let $I$ be a set and let $(Z_i)_{i\in I}$ be a family of Banach spaces. For $1\leq q<\infty$, define
\[
\ell^q\bigl(I;(Z_i)_{i\in I}\bigr)
:=
\left\{
 z=(z_i)_{i\in I}
 \ \middle|\
 z_i\in Z_i\text{ for every }i\in I,
 \quad
 \sum_{i\in I}\|z_i\|_{Z_i}^q<\infty
\right\},
\]
with the norm
\[
\|z\|_{\ell^q(I;(Z_i)_{i\in I})}
:=
\left(\sum_{i\in I}\|z_i\|_{Z_i}^q\right)^{\frac{1}{q}}.
\]
A sum over an arbitrary set is understood as the supremum of sums over finite subsets. For $q=\infty$, set
\[
\ell^\infty\bigl(I;(Z_i)_{i\in I}\bigr)
:=
\left\{
 z=(z_i)_{i\in I}
 \ \middle|\
 z_i\in Z_i\text{ for every }i\in I,
 \quad
 \sup_{i\in I}\|z_i\|_{Z_i}<\infty
\right\},
\]
with the supremum norm. If all the spaces equal a single space $Z$, write simply $\ell^q(I;Z)$. When $I=\mathbb Z$ and the family is indicated in the argument, we shall also use the abbreviation $\ell^q((Z_k)_{k\in\mathbb Z})$.
\end{definition}

\begin{lemma}[Duality of $\ell^q(I)$-sums]
\label{lem:dualidad-sumas-lq-banach}
Let $I$ be a set, let $(Z_i)_{i\in I}$ be a family of Banach spaces, and let $1\leq q<\infty$. Then
\[
\left(\ell^q\bigl(I;(Z_i)_{i\in I}\bigr)\right)'
=
\ell^{q'}\bigl(I;(Z_i')_{i\in I}\bigr)
\]
isometrically, through the pairing
\[
\langle z,\lambda\rangle
:=
\sum_{i\in I}\lambda_i(z_i).
\]
\end{lemma}

\begin{proof}
If $(\lambda_i)_{i\in I}$ belongs to the space on the right-hand side,
Hölder's inequality for sums defines a continuous functional by
\[
\Lambda((z_i)_{i\in I})=\sum_{i\in I}\lambda_i(z_i),
\qquad
\|\Lambda\|\leq\|(\|\lambda_i\|)_{i\in I}\|_{\ell^{q'}(I)}.
\]
The sum converges absolutely. Indeed, every element of
$\ell^q(I;(Z_i)_{i\in I})$, with $q<\infty$, has countable support:
each set $\{i\in I\mid\|z_i\|\geq\displaystyle\frac{1}{m}\}$ is finite, and their union
over $m\in\mathbb N$ contains every nonzero coordinate.

Conversely, let $\Lambda$ be a continuous functional. If $\iota_i$
inserts coordinate $i$, define
$\lambda_i:=\Lambda\circ\iota_i\in Z_i'$ and $a_i:=\|\lambda_i\|$.
Fix a finite set $F\subseteq I$ and $0<\varepsilon<1$.
For each $i\in F$, choose $z_i$ in the unit ball of $Z_i$ such
that $\lambda_i(z_i)$ is real and nonnegative and
\[
\lambda_i(z_i)\geq(1-\varepsilon)a_i.
\]
When $a_i=0$, take $z_i=0$; this includes $Z_i=\{0\}$.
In the complex case, multiply an approximating vector by a phase
of modulus one to obtain the stated real value.

If $1<q<\infty$ and $A_F=(\sum_{i\in F}a_i^{q'})^{\displaystyle\frac{1}{q'}}>0$,
define $c_i=\displaystyle\frac{a_i^{q'-1}}{A_F^{q'-1}}$. Then
\[
\sum_{i\in F}c_i^q=1,\qquad
\sum_{i\in F}c_i a_i=A_F.
\]
For $q=1$, choose $c_i=1$ at a coordinate maximizing $a_i$
and zero at the others; the last sum is then
$A_F=\max_{i\in F}a_i$. If $A_F=0$, there is nothing to prove.
The family with coordinates $c_i z_i$ has norm at most one, and
\[
(1-\varepsilon)A_F
\leq\Lambda\left(\sum_{i\in F}\iota_i(c_i z_i)\right)
\leq\|\Lambda\|.
\]
Let $\varepsilon\to0^+$ and take the supremum over finite sets
$F$ to obtain
\[
\|(a_i)_{i\in I}\|_{\ell^{q'}(I)}\leq\|\Lambda\|.
\]
Finally, finite-support families are dense: enumerate the nonzero
coordinates of an element and note that the tails of the sum
of their norms to the power $q$ tend to zero. On those families,
$\Lambda$ is represented by the sum of its coordinates.
By continuity, the same representation holds on the whole space.
The two norm bounds prove isometry, and the insertions $\iota_i$
prove uniqueness of every $\lambda_i$.
\end{proof}

A Banach pair $(X_0,X_1)$ is called \textit{regular} if $X_0\cap X_1$ is dense in each $X_j$. In that case, the duals $X_0'$ and $X_1'$ are viewed as subspaces of $(X_0\cap
X_1)'$ by restriction and form a compatible pair. A functional in $X_j'$ whose restriction to $X_0\cap X_1$ vanishes
is zero, by density and continuity; thus these restriction maps are injective.

\begin{lemma}[Dual of the intersection of a regular pair]
\label{lem:dual-interseccion-par-regular}
If $(X_0,X_1)$ is a regular pair, restriction gives the isometric identification
\[
(X_0\cap X_1)'
=
X_0'+X_1',
\]
where the intersection carries the maximum norm and the sum of the duals carries the quotient norm. More generally, for $t>0$,
\[
\|y\|_{(X_0\cap X_1,J(t,\cdot))'}
=
K(t^{-1},y;X_0',X_1').
\]
\end{lemma}

\begin{proof}
Consider the diagonal inclusion
\[
D\colon X_0\cap X_1\longrightarrow X_0\times X_1,
\qquad Dx=(x,x),
\]
where the product carries the norm $\displaystyle \max\{\|x_0\|_{X_0},t\|x_1\|_{X_1}\}$. Every functional $y$ on the diagonal extends without increasing its norm by the Hahn--Banach corollary~\ref{cor:extension-hahn-banach-preserva-norma}. A functional on the product has the form $(x_0,x_1)\mapsto y_0(x_0)+y_1(x_1)$ and norm $\|y_0\|_{X_0'}+t^{-1}\|y_1\|_{X_1'}$. Restricting it to the diagonal and taking the infimum over all extensions gives the formula. Regularity ensures that restriction $X_j'\to(X_0\cap X_1)'$ is injective, so the sum belongs to a compatible pair and is not merely an abstract quotient.
\end{proof}

\section{Duality for real interpolation}
\label{sec:dualidad-interpolacion-real}

\begin{theorem}[Duality of the real method]\label{teo:dualidad-interpolacion-real}\index{duality!real interpolation}
Let $(X_0,X_1)$ be a regular pair, let $0<\theta<1$, and let $1\leq q<\infty$. Then the pairing between $X_0+X_1$ and $X_0'\cap X_1'$ extends to a topological isomorphism
\[
\bigl((X_0,X_1)_{\theta,q}\bigr)'
\cong
(X_0',X_1')_{\theta,q'},
\]
where $q'$ is the conjugate exponent. For $q=\infty$, the continuous embedding of the space on the right into the dual on the left is retained, but surjectivity may fail. More precisely, if $1\leq q<\infty$ and $\Lambda_y(x):=\langle x,y\rangle$, there are constants $c_{\theta,q},C_{\theta,q}>0$ such that
\[
c_{\theta,q}\|y\|_{(X_0',X_1')_{\theta,q'}}
\leq\|\Lambda_y\|_{((X_0,X_1)_{\theta,q})'}
\leq C_{\theta,q}\|y\|_{(X_0',X_1')_{\theta,q'}}.
\]
\end{theorem}

\begin{proof}
Set $\Delta:=X_0\cap X_1$, $\Sigma:=X_0+X_1$, and $Y:=(X_0',X_1')_{\theta,q'}$. Regularity allows all elements of $Y$ to be viewed as functionals on $\Delta$.

If $a\in\Delta$ and $v\in X_0'\cap X_1'$, every decomposition $a=a_0+a_1$ satisfies
\[
|\langle a,v\rangle|
\leq\bigl(\|a_0\|_{X_0}+t^{-1}\|a_1\|_{X_1}\bigr)
      J(t,v;X_0',X_1').
\]
Thus,
\begin{equation}
|\langle a,v\rangle|
\leq K(t^{-1},a;X_0,X_1)J(t,v;X_0',X_1').
\label{eq:emparejamiento-K-J-dual}
\end{equation}
Represent $y\in Y$ by the $J$-method for the dual pair:
\[
y=\int_0^\infty v(t)\,\frac{dt}{t}
\quad\text{in }X_0'+X_1'=\Delta'.
\]
The identification of the sum space with $\Delta'$ comes from Lemma~\ref{lem:dual-interseccion-par-regular}. Evaluation at $a\in\Delta$ is continuous on this space, so \eqref{eq:emparejamiento-K-J-dual}, Hölder's inequality, and the substitution $r=t^{-1}$ give
\[
\begin{aligned}
|\langle a,y\rangle|
&\leq\int_0^\infty
 \bigl(t^\theta K(t^{-1},a)\bigr)
 \bigl(t^{-\theta}J(t,v(t))\bigr)\,\frac{dt}{t}\\
&\leq\|a\|_{\theta,q;K}
       \|t^{-\theta}J(t,v(t))\|_{L^{q'}((0,\infty),\frac{dt}{t})}.
\end{aligned}
\]
Take the infimum over representations and apply Theorem~\ref{teo:equivalencia-metodos-K-J}. If $q<\infty$, density of $\Delta$ from Lemma~\ref{lem:caracterizacion-discreta-metodo-J} allows $y$ to extend uniquely to a continuous functional $\Lambda_y$ on $(X_0,X_1)_{\theta,q}$, with the asserted bound.

Conversely, let $\Lambda\in((X_0,X_1)_{\theta,q})'$ and assume $q<\infty$. Define
\[
Z_k:=(\Delta,2^{-k\theta}J(2^k,\cdot)),\qquad
\mathcal J_{\theta,q}:=\ell^q\bigl(\mathbb Z;(Z_k)_k\bigr).
\]
Each $Z_k$ is Banach by Proposition~\ref{prop:interseccion-suma-banach}. Lemma~\ref{lem:caracterizacion-discreta-metodo-J} shows that $\displaystyle S((a_k)_{k\in\mathbb Z}):=\displaystyle\sum_{k\in \mathbb Z} a_k$ is continuous and surjective from $\mathcal J_{\theta,q}$ onto the interpolated space. Apply Lemma~\ref{lem:dualidad-sumas-lq-banach} directly to the functional $\Lambda\circ S$. We obtain functionals $\lambda_k\in\Delta'$ such that
\[
\left\|\left(
2^{k\theta}\|\lambda_k\|_{(\Delta,J(2^k,\cdot))'}
\right)_k\right\|_{\ell^{q'}(\mathbb Z)}\leq C_{\theta,q}\|\Lambda\|.
\]
Evaluation on a sequence whose only nonzero coordinate is $a\in\Delta$ gives $\lambda_k(a)=\Lambda(a)$. Thus all $\lambda_k$ agree with $y:=\Lambda|_\Delta$. Lemma~\ref{lem:dual-interseccion-par-regular} transforms the estimate into
\[
\left\|\left(2^{k\theta}K(2^{-k},y;X_0',X_1')\right)_k
\right\|_{\ell^{q'}(\mathbb Z)}\leq C_{\theta,q}\|\Lambda\|.
\]
Lemma~\ref{lem:discretizacion-funcional-K}, with the change of index $k\mapsto-k$, implies $y\in Y$ and gives the other norm bound. Since $\Delta$ is dense, $\Lambda=\Lambda_y$ on the whole space.

For $q=\infty$, the pairing outside $\Delta$ requires clarification. Let $x\in(X_0,X_1)_{\theta,\infty}$ and take the partial sums $x_N$ of a nearly optimal discrete $J$-representation. Then $x_N\in\Delta$, $x_N\to x$ in $\Sigma$, and
\[
\sup_{N\in\mathbb N}\|x_N\|_{\theta,\infty;K}
\leq C_\theta\|x\|_{\theta,\infty;K}.
\]
If $y\in X_0'\cap X_1'$, its action is continuous on $\Sigma$. The estimate already obtained on $\Delta$, followed by the limit in $N$, gives
\[
|\langle x,y\rangle|
\leq C_\theta\|x\|_{\theta,\infty;K}
                  \|y\|_{(X_0',X_1')_{\theta,1}}.
\]
Lemma~\ref{lem:caracterizacion-discreta-metodo-J}, applied to the dual pair with second index $1$, extends this pairing to all of $y\in(X_0',X_1')_{\theta,1}$. It is bilinear, continuous, and injective in $y$, since on $\Delta$ it agrees with the original action. The preceding proof of surjectivity requires duality of $\ell^q(\mathbb Z)$ with $q<\infty$ and does not apply at the endpoint $q=\infty$.
\end{proof}

\section{The complex method}
\label{sec:metodo-complejo-interpolacion}

In this section, the Banach spaces are complex. A function $F\colon U\subseteq\mathbb C\longrightarrow X$ is called holomorphic if the limit
\[
F'(z):=\lim_{h\to0}\frac{F(z+h)-F(z)}{h}
\]
exists in the norm of $X$ for each $z\in U$.

To apply the method to a real space $X$, use its complexification $X_{\mathbb C}:=X+iX$ with the norm
\[
\|x+iy\|_{X_{\mathbb C}}
:=\sup_{t\in\mathbb R}\|x\cos t-y\sin t\|_X.
\]
Multiplication by a complex scalar rotates the pair $(x,y)$ and multiplies it by the modulus of the scalar, so this expression is a complex norm. Moreover,
\[
\max\{\|x\|_X,\|y\|_X\}
\leq\|x+iy\|_{X_{\mathbb C}}\leq\|x\|_X+\|y\|_X,
\]
which proves completeness. The inclusions of a real pair are complexified within a common ambient space. Define its real interpolated space as the fixed-point space of conjugation on $[X_{0,\mathbb C},X_{1,\mathbb C}]_\theta$. Conjugation preserves the Calderón norm upon replacing $F(z)$ by $\overline{F(\bar z)}$. If $T$ is real and continuous, its complex extension has the same norm: the inequality $\|T_{\mathbb C}(x+iy)\|\leq\|T\|\|x+iy\|$ follows for each term in the supremum, and the reverse bound follows on $X$. Thus the operator theorems also apply to real spaces. The isometric equalities below are formulated for complex spaces with the indicated norms. Passage to their real parts retains the identifications with equivalent norms.

We shall denote by
\[
\mathcal T:=\{z\in\mathbb C\mid0<\operatorname{Re}z<1\},
\qquad
\overline{\mathcal T}
=\{z\in\mathbb C\mid0\leq\operatorname{Re}z\leq1\}
\]
the open strip and its closure, respectively.

Before defining the Calderón space, we establish the complex analysis tools to be used in both the scalar and vector-valued cases.

\begin{proposition}[Cauchy's integral formula with values in a Banach space]
\label{prop:formula-integral-cauchy-banach}
Let $U\subseteq\mathbb C$ be open, let $X$ be a complex Banach space, and let $F\colon U\longrightarrow X$ be holomorphic. If $\overline B_{\mathrm{euc}}(z_0,r)\subseteq U$, then for $z\in B_{\mathrm{euc}}(z_0,r)$ and $m\in\mathbb N_0$,
\begin{equation}
F^{(m)}(z)
=
\frac{m!}{2\pi i}
\int_{|\zeta-z_0|=r}
\frac{F(\zeta)}{(\zeta-z)^{m+1}}\,d\zeta,
\label{eq:formula-cauchy-banach}
\end{equation}
where the circle is traversed with positive orientation and the integral is a Bochner integral. In particular, $F$ admits a Taylor series around each point, convergent in the norm of $X$. If $X=\mathbb C$ and two holomorphic functions agree on a set with an interior accumulation point, they agree on every connected component of $U$ containing that point.
\end{proposition}

\begin{proof}
Begin with the scalar case and Goursat's argument. If the integral of a holomorphic function over the boundary of a triangle were nonzero, subdividing the triangle into four congruent triangles would allow us to choose one whose integral has modulus at least one quarter of the original. Iteration produces nested triangles with diameters tending to zero and a common point $a$. For each $\varepsilon>0$, differentiability of $f$ at $a$ provides $\delta>0$ such that
\[
 |f(\zeta)-f(a)-f'(a)(\zeta-a)|
 \leq\varepsilon|\zeta-a|
 \qquad\text{if }|\zeta-a|<\delta.
\]
The integrals of the two affine terms over a triangular boundary vanish. If $d_k:=\operatorname{diam}(T_k)<\delta$, every point of $T_k$ is at distance at most $d_k$ from $a$, and the length of $\partial T_k$ does not exceed $3d_k$; consequently,
\[
 \left|\int_{\partial T_k}f(\zeta)\,d\zeta\right|
 \leq3\varepsilon d_k^2.
\]
On the other hand, the triangles are congruent after rescaling by $2^{-k}$, so $d_k=2^{-k}d_0$, whereas the inductive choice implies
\[
 \left|\int_{\partial T_k}f(\zeta)\,d\zeta\right|
 \geq4^{-k}
 \left|\int_{\partial T_0}f(\zeta)\,d\zeta\right|.
\]
If the last integral were nonzero, division by $d_k^2$ and the choice $\displaystyle 3\varepsilon<d_0^{-2}
|\int_{\partial T_0}f(\zeta)\,d\zeta|$ would give a contradiction. Thus the integral over every triangular boundary is zero. Summing over a triangulation cancels the integrals along interior edges. This gives the result for polygonal regions, including those with holes. For the circles we shall use, approximate by inscribed polygons; uniform continuity of the integrand on a compact neighborhood and convergence of these parametrizations allow passage to the limit in the integrals.

Apply this fact to $\zeta\mapsto\displaystyle\frac{f(\zeta)}{\zeta-z}$ on the region between $|\zeta-z_0|=r$ and a small circle centered at $z$. Shrinking the latter circle yields
\[
f(z)
=
\frac{1}{2\pi i}
\int_{|\zeta-z_0|=r}\frac{f(\zeta)}{\zeta-z}\,d\zeta.
\]
Differentiation of the kernel under the integral gives the formula for $f^{(m)}$.

Return to the vector-valued case. Continuity of $F$ on the circle ensures existence of the Bochner integral. For every $\lambda\in X'$, the function $\lambda\circ F$ is holomorphic, and compatibility of the integral with continuous linear operators, Proposition~\ref{prop:b5-integral-bochner-operadores-lineales}, allows application of the scalar formula to $m=0$. The Hahn--Banach theorem \ref{teo: hahn--banach, forma analitica} shows that $X'$ separates points and gives the vector-valued identity for $F(z)$. On each concentric disk of radius less than $r$, the difference quotients of $(\zeta-z)^{-1}$ converge uniformly on the circle to $(\zeta-z)^{-2}$. The Bochner integral inequality allows differentiation of the identity in the norm of $X$. Repetition gives all derivatives and the formula for each $m$.

If $|z-z_0|<\rho<r$, the geometric series for the kernel converges uniformly for $|z-z_0|\leq\rho$; we may integrate term by term and obtain the Taylor series of $F$, converging uniformly in norm on that disk. In the scalar case, the first nonzero coefficient in this series shows that the zeros of a holomorphic function that is not identically zero are isolated. This proves the final assertion.
\end{proof}

\begin{lemma}[Maximum principle and holomorphic limits]
\label{lem:maximo-limites-holomorfos}
The following assertions hold.
\begin{enumerate}[label=(\alph*)]
\item If $Q\subseteq\mathbb C$ is a compact rectangle and $f\colon Q\longrightarrow\mathbb C$ is continuous and holomorphic in the interior, then $\displaystyle\max_{z\in Q}|f(z)|=\displaystyle\max_{z\in\partial Q}|f(z)|$. Moreover, if $|f|$ attains its maximum at an interior point, then $f$ is constant.
\item A uniform limit on compact subsets of scalar holomorphic functions is holomorphic.
\item If $X$ is a Banach space, a uniform limit on compact subsets of holomorphic functions with values in $X$ is holomorphic.
\end{enumerate}
\end{lemma}

\begin{proof}
Proposition~\ref{prop:formula-integral-cauchy-banach}, with $m=0$, gives the mean value property. If $|f|$ attains an interior maximum, equality in the triangle inequality for that formula forces $f$ to be constant on a neighborhood; the final assertion of the same proposition propagates constancy throughout the connected interior of the rectangle. This proves (a).

For (b) and (c), let $(F_j)$ be a sequence converging uniformly on compact subsets to $F$, with values in $\mathbb C$ or $X$, respectively. Each $F_j$ satisfies Cauchy's formula from Proposition~\ref{prop:formula-integral-cauchy-banach}. Uniform convergence on the circle allows passage to the limit in the Bochner integral and gives
\[
F(z)=\frac{1}{2\pi i}\int_{|\zeta-z_0|=r}
\frac{F(\zeta)}{\zeta-z}\,d\zeta
\]
on every disk whose closure is contained in the common domain. On a smaller concentric disk, the kernel and its derivative with respect to $z$ are uniformly bounded on the circle. The fundamental theorem of calculus~\ref{teo:b5-fundamental-calculo-riemann-banach} then allows differentiation under the integral and proves holomorphicity of $F$.
\end{proof}

\begin{lemma}[Hadamard's three-lines lemma]\label{lem:tres-lineas-hadamard}\index{three-lines lemma}
Let $F\colon \overline{\mathcal T}\longrightarrow\mathbb C$ be continuous and bounded, holomorphic on $\mathcal T$, and suppose that $|F(j+it)|\leq M_j$ for $t\in\mathbb R$ and $j=0,1$. Then
\[
|F(\theta)|
\leq
M_0^{1-\theta}M_1^\theta,
\qquad 0<\theta<1.
\]
\end{lemma}

\begin{proof}
We may assume $M_0M_1>0$. For $\varepsilon>0$, consider
\[
G_\varepsilon(z)
:=
F(z)M_0^{z-1}M_1^{-z}
e^{\varepsilon(z-\theta)^2}.
\]
On the vertical sides of a rectangle $0\leq\operatorname{Re}z\leq1$, $|\operatorname{Im}z|\leq R$, its modulus is bounded by $\displaystyle e^{\varepsilon\displaystyle\max\{\theta^2,(1-\theta)^2\}}$; on the horizontal sides, it tends to zero as $R\to\infty$. The maximum principle from Lemma~\ref{lem:maximo-limites-holomorfos}, applied to the rectangle and followed by the limit $R\to\infty$, gives $\displaystyle |G_\varepsilon(\theta)|\leq
e^{\varepsilon\displaystyle\max\{\theta^2,(1-\theta)^2\}}$. Finally, let $\varepsilon\to0^+$. If $M_j=0$, apply the preceding argument with $M_j+\delta$ and let $\delta\to0^+$.
\end{proof}

The two lines bounding \(\mathcal T\) represent the endpoint spaces: on \(it\) we measure in \(X_0\), and on \(1+it\) in \(X_1\), for \(t\in\mathbb R\). A holomorphic function joins these boundary values. Evaluation at \(z=\theta\) produces the intermediate element, and the three-lines lemma controls its size.

\begin{definition}[Calderón space]\label{def:espacio-calderon-interpolacion-compleja}\index{complex interpolation}
\glsadd{interpolacion-compleja}
Let $(X_0,X_1)$ be a compatible pair of complex Banach spaces. Denote by $\mathscr
F(X_0,X_1)$ the space of functions $F\colon \overline{\mathcal T}\longrightarrow X_0+X_1$ that are continuous and bounded, holomorphic on $\mathcal T$, satisfy $F(j+it)\in X_j$, have restrictions $t\mapsto F(j+it)$ continuous in the norm of $X_j$ with $F(j+it)\to0$ in $X_j$ as $|t|\to\infty$, and have norm
\[
\|F\|_{\mathscr F}
:=
\max_{j\in\{0,1\}}\sup_{t\in\mathbb R}\|F(j+it)\|_{X_j}<\infty.
\]
For $0<\theta<1$, define
\[
[X_0,X_1]_\theta
:=
\{F(\theta)\mid F\in\mathscr F(X_0,X_1)\},
\qquad
\|x\|_{[X_0,X_1]_\theta}
:=
\inf_{\substack{F\in\mathscr F(X_0,X_1)\\F(\theta)=x}}\|F\|_{\mathscr F}.
\]
\end{definition}

\begin{theorem}[Complex interpolation of operators]\label{teo:interpolacion-compleja-operadores}\index{interpolation of operators!complex method}
The space $[X_0,X_1]_\theta$ is Banach. If $T\colon X_0+X_1\longrightarrow Y_0+Y_1$ is linear and $T\colon X_j\longrightarrow Y_j$ has norm $M_j$, then
\[
T\colon [X_0,X_1]_\theta\longrightarrow[Y_0,Y_1]_\theta
\]
is continuous and $\|T\|\leq M_0^{1-\theta}M_1^\theta$.
\end{theorem}

\begin{proof}
Write $\Sigma:=X_0+X_1$. For $\lambda\in\Sigma'$ of norm at most one, the restrictions of $\lambda$ to $X_j$ have norm at most one. Lemma~\ref{lem:tres-lineas-hadamard}, applied to $z\mapsto\lambda(F(z+it))$, and Hahn--Banach give
\begin{equation}
\sup_{z\in\overline{\mathcal T}}\|F(z)\|_\Sigma
\leq\|F\|_{\mathscr F}.
\label{eq:control-interior-calderon}
\end{equation}
Thus a Cauchy sequence in $\mathscr F$ converges uniformly in $\Sigma$ on the closed strip and uniformly in $X_j$ on each boundary line. The limit is continuous in the indicated norms, its boundary values tend to zero, and Lemma~\ref{lem:maximo-limites-holomorfos} preserves holomorphicity. Hence $\mathscr F$ is Banach. Evaluation $F\mapsto F(\theta)$ is continuous by \eqref{eq:control-interior-calderon}, so its kernel is closed. The space $[X_0,X_1]_\theta$, with its defining norm, is the quotient by that kernel and is Banach. In particular,
\[
\|x\|_\Sigma\leq\|x\|_{[X_0,X_1]_\theta}.
\]
If $x\in X_0\cap X_1$, the representative $F(z)=e^{\varepsilon(z-\theta)^2}x$ also gives the continuous embedding $X_0\cap X_1\hookrightarrow[X_0,X_1]_\theta$.

The operator $T$ is continuous between the sum spaces: apply its endpoint bounds to a decomposition and take the infimum. If $M_0M_1>0$ and $F\in\mathscr F(X_0,X_1)$, set
\[
G(z):=M_0^{z-1}M_1^{-z}T(F(z)).
\]
This function belongs to $\mathscr F(Y_0,Y_1)$ and satisfies $\|G\|_{\mathscr F}\leq\|F\|_{\mathscr F}$. Since $G(\theta)=M_0^{\theta-1}M_1^{-\theta}T(F(\theta))$, the infimum over $F$ gives the required bound. If some $M_j$ is zero, use the bounds $M_j+\varepsilon$ and let $\varepsilon\to0^+$.
\end{proof}

\begin{lemma}[Dual criterion for the Lebesgue norm]
\label{lem:criterio-dual-Lp-interpolacion}
Let $(A,\mu)$ be a $\sigma$-finite measure space and let $h$ be a measurable function finite almost everywhere. For $1\leq s\leq\infty$,
\[
\|h\|_{L^s(A)}
=\sup\left\{\left|\int_A h g\,d\mu\right|
\;\middle|\;
\begin{array}{l}
g\text{ simple, supported on a set of finite measure},\\
hg\in L^1(A),\quad\|g\|_{L^{s'}(A)}\leq1
\end{array}\right\},
\]
allowing both sides to be infinite.
\end{lemma}

\begin{proof}
Hölder's inequality proves one bound if $h\in L^s(A)$; if the norm of $h$ is infinite, this bound holds automatically. For the reverse bound, take an increasing sequence of sets of finite measure whose union is $A$ and intersect them with $\{|h|\leq N\}$. On each resulting set, when $1<s<\infty$, use the function $\overline{\operatorname{sgn}h}|h|^{s-1}$ normalized in $L^{s'}(A)$. It is bounded and can be approximated uniformly by simple functions; if necessary, divide each approximant by the maximum of one and its norm. The supremum is at least the $L^s(A)$ norm of $h$ restricted to that set. For $s=1$, use $\overline{\operatorname{sgn}h}$, of norm at most one in $L^\infty(A)$. Monotone convergence gives the conclusion in both cases. If $s=\infty$, for each $0<a<\|h\|_\infty$ choose a set of positive finite measure on which $a<|h|\leq N$, and use its normalized indicator function multiplied by $\overline{\operatorname{sgn}h}$. The same simple approximations show that the supremum is at least $a$.
\end{proof}

\begin{theorem}[Riesz--Thorin]\label{teo:riesz-thorin-interpolacion}\index{Riesz--Thorin theorem}
Let $(A,\mu)$ and $(B,\nu)$ be $\sigma$-finite measure spaces, let $1\leq p_0,p_1,q_0,q_1\leq\infty$, and let
\[
T:L^{p_0}(A)+L^{p_1}(A)
\longrightarrow
L^{q_0}(B)+L^{q_1}(B)
\]
be a linear operator. Suppose that
\[
\|Tf\|_{L^{q_j}(B)}
\leq
M_j\|f\|_{L^{p_j}(A)},
\qquad j=0,1.
\]
If $0<\theta<1$ and
\[
\frac{1}{p_\theta}
=
\frac{1-\theta}{p_0}+\frac{\theta}{p_1},
\qquad
\frac{1}{q_\theta}
=
\frac{1-\theta}{q_0}+\frac{\theta}{q_1},
\]
then the restriction of $T$ to $L^{p_\theta}(A)$ takes values in $L^{q_\theta}(B)$ and
\[
\|Tf\|_{L^{q_\theta}(B)}
\leq
M_0^{1-\theta}M_1^\theta\|f\|_{L^{p_\theta}(A)}.
\]
\end{theorem}

\begin{proof}
First assume $p_\theta<\infty$ and $q_\theta>1$. By homogeneity, it suffices to consider simple functions
\[
f=\displaystyle\sum_{k=1}^{N}a_k\mathbf1_{E_k},
\qquad
g=\displaystyle\sum_{\ell=1}^{L}b_\ell\mathbf1_{F_\ell},
\]
with disjoint sets of finite measure and $\|f\|_{L^{p_\theta}(A)}=
\|g\|_{L^{q_\theta'}(B)}=1$. For $a\neq0$, write $\omega(a):=\displaystyle\frac{a}{|a|}$ and set
\[
\frac{1}{p(z)}:=\frac{1-z}{p_0}+\frac{z}{p_1},
\qquad
\frac{1}{q'(z)}:=\frac{1-z}{q_0'}+\frac{z}{q_1'}.
\]
With the convention that zero coefficients are omitted, define
\[
f_z:=\displaystyle\sum_{k=1}^{N}
|a_k|^{\frac{p_\theta}{p(z)}}\omega(a_k)\mathbf1_{E_k},
\qquad
g_z:=\displaystyle\sum_{\ell=1}^{L}
|b_\ell|^{\frac{q_\theta'}{q'(z)}}\omega(b_\ell)\mathbf1_{F_\ell}.
\]
Then $f_\theta=f$, $g_\theta=g$, and, for $j\in\{0,1\}$ and $t\in\mathbb R$,
\[
\|f_{j+it}\|_{L^{p_j}(A)}=1,
\qquad
\|g_{j+it}\|_{L^{q_j'}(B)}=1.
\]
If $p_j=\infty$ or $q_j'=\infty$, the equality is interpreted with $\displaystyle\frac{1}{\infty}=0$ and verified by taking the essential maximum. The function
\[
\Phi(z):=
\int_B(Tf_z)(y)g_z(y)\,d\nu(y)
\]
is holomorphic and bounded on the closed strip; expanding the two sums reduces it to a finite sum of exponential functions of $z$. The boundary hypothesis and Hölder's inequality from Proposition~\ref{desigualdad de holder} give $|\Phi(j+it)|\leq M_j$. Lemma~\ref{lem:tres-lineas-hadamard} yields
\[
\left|\int_B(Tf)g\,d\nu\right|
\leq
M_0^{1-\theta}M_1^\theta.
\]
Undoing the normalization and using the criterion in Lemma~\ref{lem:criterio-dual-Lp-interpolacion} gives the estimate for $1<q_\theta<\infty$. If $q_\theta=1$, necessarily $q_0=q_1=1$; take $g_z=g\in L^\infty(B)$, and the same argument works. If $q_\theta=\infty$, necessarily $q_0=q_1=\infty$; obtain the bound by applying the argument to functions $g\in L^1(B)$ and using duality again. If $p_\theta=\infty$, then $p_0=p_1=\infty$, and take $f_z=f$.

Finally, Proposition~\ref{prop:densidad-funciones-simples-Lp} allows passage from simple functions to every function in $L^{p_\theta}(A)$ when $p_\theta<\infty$. This limit preserves the original operator: $L^{p_\theta}(A)$ embeds continuously into $L^{p_0}(A)+L^{p_1}(A)$. To see this, if $p_0<p_1$ and $f\neq0$, split $f$ over the sets $\{|f|>\|f\|_{p_\theta}\}$ and $\{|f|\leq\|f\|_{p_\theta}\}$; the two parts have norms in $L^{p_0}(A)$ and $L^{p_1}(A)$ at most $\|f\|_{p_\theta}$. If $p_1<p_0$, interchange the roles, and if $p_0=p_1$ the embedding is immediate. The operator $T$ is continuous between the sum spaces by its endpoint bounds, and both limits agree locally in measure. When $p_\theta=\infty$, the argument with $f_z=f$ applies to each $f\in L^\infty(A)$ without approximation in that norm. This completes the proof.
\end{proof}

\begin{theorem}[Hausdorff--Young inequality]\label{teo:hausdorff-young-fourier}\index{Hausdorff--Young inequality}
Let $1\leq p\leq2$ and let $p'$ be its conjugate exponent. The Fourier transform extends uniquely to a continuous operator
\[
\mathcal F\colon L^p(\mathbb R^n)\longrightarrow L^{p'}(\mathbb R^n)
\]
and satisfies
\[
\|\widehat f\|_{L^{p'}(\mathbb R^n)}
\leq
(2\pi)^{-n(\frac{1}{p}-\frac{1}{2})}\|f\|_{L^p(\mathbb R^n)}.
\]
\end{theorem}

\begin{proof}
The norm of the operator from $L^1(\mathbb R^n)$ to $L^\infty(\mathbb R^n)$ is at most $(2\pi)^{-\displaystyle\frac{n}{2}}$ by Proposition~\ref{prop:fourier-L1-continuidad}; on $L^2(\mathbb R^n)$ it is one by Plancherel's theorem~\ref{teo:plancherel-fourier}. The Riesz--Thorin theorem~\ref{teo:riesz-thorin-interpolacion}, applied initially to simple functions of finite support, gives the result for $1<p<2$ with the indicated constant. Density of that subspace in $L^p(\mathbb R^n)$, given by Proposition~\ref{prop:densidad-funciones-simples-Lp}, provides the unique extension.
\end{proof}

The interval $1\leq p\leq2$ in Hausdorff--Young is essential. For $p>2$, Proposition~\ref{prop:Lp-distribucion-temperada-fourier} still defines the transform in $\mathcal S'(\mathbb R^n)$, but in general there is no continuous operator $L^p(\mathbb R^n)\to L^{p'}(\mathbb R^n)$. For $p=\infty$, the simplest instance already appears: the constant function one satisfies
\[
\mathcal F1
=
(2\pi)^{\frac{n}{2}}\delta_0
\]
in $\mathcal S'(\mathbb R^n)$, which is not a function. For $2<p<\infty$, failure of a bound $L^p(\mathbb R^n)\to L^{p'}(\mathbb R^n)$ can be seen through sums of modulations. If $\psi\in\mathcal S(\mathbb R^n)$, $\widehat\psi$ is supported in a sufficiently small ball, and $L>0$ separates its translates, consider
\[
f_N(x)
:=
\psi(x)\displaystyle\sum_{k=1}^N\varepsilon_ke^{ikLx_1},
\qquad \varepsilon_k\in\{-1,1\}.
\]
Theorem~\ref{teo:khintchine}, integrated in $x$ by Theorem~\ref{teo:tonelli}, allows the signs to be chosen so that $\|f_N\|_{L^p(\mathbb R^n)}\leq C_pN^{\displaystyle\frac{1}{2}}\|\psi\|_{L^p(\mathbb R^n)}$. Disjointness of the supports of the frequency blocks gives $\|\widehat f_N\|_{L^{p'}(\mathbb R^n)}=N^{\displaystyle\frac{1}{p'}}\|\widehat\psi\|_{L^{p'}(\mathbb
R^n)}$. Since $\displaystyle\frac{1}{p'}>\displaystyle\frac{1}{2}$, a uniform bound would be impossible.

\begin{corollary}[Complex interpolation of Lebesgue spaces]
\label{cor:interpolacion-compleja-Lp}
Let $(A,\mu)$ be a $\sigma$-finite measure space, let $1\leq p_0,p_1\leq\infty$, and let $0<\theta<1$. If
\[
\frac{1}{p}=\frac{1-\theta}{p_0}+\frac{\theta}{p_1},
\]
then
\[
[L^{p_0}(A),L^{p_1}(A)]_\theta=L^p(A)
\]
with equality of norms under the usual normalization of the complex method.
\end{corollary}

\begin{proof}
Let $F\in\mathscr F(L^{p_0}(A),L^{p_1}(A))$. For a simple function $g$ supported on a set of finite measure, with $\|g\|_{L^{p'}(A)}=1$, construct as in the proof of Theorem~\ref{teo:riesz-thorin-interpolacion} an analytic family $g_z$ such that $g_\theta=g$ and $\|g_{j+it}\|_{L^{p_j'}(A)}=1$. The three-lines lemma \ref{lem:tres-lineas-hadamard}, applied to
\[
z\longmapsto\int_A F(z,x)g_z(x)\,d\mu(x)
\]
gives $\left|\displaystyle\int_A F(\theta,x)g(x)\,d\mu(x)\right|
\leq\|F\|_{\mathscr F}$. The dual criterion in Lemma~\ref{lem:criterio-dual-Lp-interpolacion} yields
\[
\|F(\theta)\|_{L^p(A)}\leq\|F\|_{\mathscr F},
\]
and taking the infimum gives $[L^{p_0}(A),L^{p_1}(A)]_\theta\hookrightarrow L^p(A)$.

For the reverse inclusion, if $p<\infty$, begin with a nonzero simple function $f$ supported on a set of finite measure and consider
\[
F(z,x)
:=
e^{\varepsilon(z-\theta)^2}
\|f\|_{L^p(A)}^{1-\frac{p}{p(z)}}
|f(x)|^{\frac{p}{p(z)}}\operatorname{sgn}f(x),
\qquad
\frac{1}{p(z)}=\frac{1-z}{p_0}+\frac{z}{p_1}.
\]
Then $F(\theta)=f$ and, on the boundary $j+it$, the modulus of the power of $f$ is $|f|^{\displaystyle\frac{p}{p_j}}$; hence
\[
\|F(j+it)\|_{L^{p_j}(A)}
\leq
e^{\varepsilon(j-\theta)^2}\|f\|_{L^p(A)}.
\]
The Gaussian factor also ensures decay as $|t|\to\infty$. Letting $\varepsilon\to0^+$ gives
\[
\|f\|_{[L^{p_0}(A),L^{p_1}(A)]_\theta}
\leq\|f\|_{L^p(A)}.
\]
Proposition~\ref{prop:densidad-funciones-simples-Lp} and completeness from Theorem~\ref{teo:interpolacion-compleja-operadores} allow passage to every $f\in L^p(A)$ when $p<\infty$. If an endpoint exponent is infinite, use the convention $\displaystyle\frac{1}{\infty}=0$ in the same formula; if $p=\infty$, necessarily $p_0=p_1=\infty$: the function $z\mapsto
e^{\varepsilon(z-\theta)^2}f$ gives one bound, and the three-lines lemma gives the other.
\end{proof}

\subsection{Translations, regularization, and boundary estimates}

To interpolate families of spaces, we shall need more than a bound by the maximum of the boundary norms. We first justify the approximations to be used, then obtain an estimate averaging those norms.

\begin{lemma}[Translations and Gaussian regularization]
\label{lem:regularizacion-gaussiana-calderon}
Let $F\in\mathscr F(X_0,X_1)$ and $0<\alpha<1$. Then
\[
t\longmapsto F(\alpha+it)
\]
is continuous with values in $[X_0,X_1]_\alpha$, tends to zero as $|t|\to\infty$, and satisfies
\[
\sup_{t\in\mathbb R}\|F(\alpha+it)\|_{[X_0,X_1]_\alpha}\leq\|F\|_{\mathscr F}.
\]
Moreover, there are functions $F_\varepsilon\in\mathscr F(X_0,X_1)$, entire with values in $X_0\cap X_1$, such that $F_\varepsilon\to F$ in $\mathscr F$. In particular, $X_0\cap X_1$ is dense in $[X_0,X_1]_\alpha$.
\end{lemma}

\begin{proof}
The function $z\mapsto F(z+it)$ belongs to $\mathscr F$ with the same norm as $F$, giving the bound. The boundary functions are continuous and tend to zero, so they are uniformly continuous. Hence
\[
\|F(\,\cdot+ih)-F\|_{\mathscr F}\longrightarrow0\qquad(h\to0),
\]
and evaluation at $\alpha+it$ proves the asserted continuity. To justify the limit at infinity, fix $d>0$ and use the representative $z\mapsto e^{d(z-\alpha)^2}F(z+it)$. We obtain
\[
\|F(\alpha+it)\|_{[X_0,X_1]_\alpha}
\leq e^{d\max\{\alpha^2,(1-\alpha)^2\}}
\max_{j\in\{0,1\}}\sup_{s\in\mathbb R} e^{-ds^2}\|F(j+i(t+s))\|_{X_j}.
\]
On $|s|>R$, the Gaussian factor makes the right-hand side small uniformly in $t$; on $|s|\leq R$, the boundary values tend uniformly to zero as $|t|\to\infty$.

For $\varepsilon>0$, define, initially on the strip,
\[
F_\varepsilon(z):=\frac{1}{\sqrt{\pi\varepsilon}}
\int_{\mathbb R}e^{-\frac{t^2}{\varepsilon}}F(z+it)\,dt.
\]
The integrals are Bochner integrals. Cauchy's theorem used in Proposition~\ref{prop:formula-integral-cauchy-banach}, applied to the integrand over rectangles, allows the integration line to be shifted and gives
\begin{equation}
F_\varepsilon(z)=\frac{1}{\sqrt{\pi\varepsilon}}
\int_{\mathbb R}e^{\frac{(z-j-it)^2}{\varepsilon}}F(j+it)\,dt,
\qquad j=0,1.
\label{eq:regularizacion-calderon-dos-fronteras}
\end{equation}
Indeed, on the horizontal sides, the modulus of the kernel is bounded by a constant times $e^{-\displaystyle\frac{(R-|\operatorname{Im}z|)^2}{\varepsilon}}$; their integrals tend to zero. First shift interior lines, then pass to the boundary by continuity and dominated convergence.

The right-hand side of \eqref{eq:regularizacion-calderon-dos-fronteras} defines an entire function with values in $X_j$: on every compact subset of $z$, all derivatives of the kernel admit a bound integrable in $t$. The two expressions agree on the strip and, by the identity principle, throughout $\mathbb C$ with values in $X_0+X_1$. They thus define an entire function in $X_0\cap X_1$. On the line $j+i\mathbb R$, they are Gaussian convolutions of the boundary function in $X_j$, whence
\[
\|F_\varepsilon\|_{\mathscr F}\leq\|F\|_{\mathscr F},\qquad
\|F_\varepsilon-F\|_{\mathscr F}\longrightarrow0.
\]
For the second assertion, split the integral into $|t|<\delta$, where uniform continuity is used, and $|t|\geq\delta$, where the Gaussian mass tends to zero. Evaluation at $\alpha$ proves density.
\end{proof}

\begin{lemma}[Integral estimate for the interpolated norm]
\label{lem:poisson-norma-interpolada}
For $0<\alpha<1$, set
\[
P_0(\alpha,t):=\frac{\sin(\pi\alpha)}{2(\cosh(\pi t)-\cos(\pi\alpha))},\qquad
P_1(\alpha,t):=\frac{\sin(\pi\alpha)}{2(\cosh(\pi t)+\cos(\pi\alpha))}.
\]
These kernels are positive and integrable, with masses $m_0:=1-\alpha$ and $m_1:=\alpha$, respectively. If $F\in\mathscr F(X_0,X_1)$, then
\begin{equation}
\|F(\alpha)\|_{[X_0,X_1]_\alpha}
\leq\exp\left(\sum_{j=0}^1\int_{\mathbb R}
 P_j(\alpha,t)\log\|F(j+it)\|_{X_j}\,dt\right).
\label{eq:poisson-log-norma-interpolada}
\end{equation}
We use $\log0=-\infty$. In particular, if $1\leq r_0,r_1\leq\infty$,
\begin{equation}
\|F(\alpha)\|_{[X_0,X_1]_\alpha}
\leq B_0^{1-\alpha}B_1^\alpha,
\label{eq:poisson-potencias-norma-interpolada}
\end{equation}
where
\[
B_j:=\left(\frac{1}{m_j}\int_{\mathbb R}
 P_j(\alpha,t)\|F(j+it)\|_{X_j}^{r_j}\,dt\right)^{\frac{1}{r_j}}
\]
if $r_j<\infty$, and $\displaystyle B_j:=\displaystyle\sup_{t\in\mathbb R}\|F(j+it)\|_{X_j}$ if $r_j=\infty$. We also have
\begin{equation}
\|F(\alpha)\|_{[X_0,X_1]_\alpha}
\leq\sum_{j=0}^1\int_{\mathbb R}
P_j(\alpha,t)\|F(j+it)\|_{X_j}\,dt.
\label{eq:poisson-lineal-norma-interpolada}
\end{equation}
The same estimates hold for a continuous bounded function on the strip, holomorphic in its interior, whose boundary restrictions are continuous and bounded in $X_j$, even if they do not tend to zero.
\end{lemma}

\begin{proof}
We first describe the harmonic function used in the proof. For smooth real functions $\varphi_0,\varphi_1$ constant outside a compact interval, set
\[
u(x,y):=\sum_{j=0}^1\int_{\mathbb R}
 P_j(x,y-t)\varphi_j(t)\,dt,\qquad 0<x<1.
\]
The kernels are obtained through the substitution $w=e^{\pi iz}$ in the Poisson kernel for the upper half-plane
\[
\frac{\operatorname{Im}w}{\pi|s-w|^2}\,ds.
\]
For $s>0$, write $s=e^{-\pi t}$, and for $s<0$, $s=-e^{-\pi t}$; these substitutions give exactly $P_0$ and $P_1$. The half-plane kernel is the imaginary part of $\pi^{-1}(s-w)^{-1}$ and is therefore harmonic. Its integral in $s$ is one, and its mass outside a neighborhood of the boundary point tends to zero. This proves harmonicity of $u$ and its boundary values $u(j,y)=\varphi_j(y)$. Integrating the kernel over the positive and negative half-lines, using an arctangent antiderivative, gives the masses $1-x$ and $x$. In particular, $u$ lies between the minimum and maximum of the boundary data. For smooth data as chosen, $u$ and its first derivatives extend continuously to the boundary lines. This can be checked in the half-plane kernel: for the tangential derivative, integrate by parts and apply the kernel to the differentiated datum. For the normal derivative, subtract the first-order Taylor polynomial of the datum near the boundary point. The linear term has zero integral over symmetric intervals; the remainder is of second order and provides an integrable bound for the differentiated kernel. Outside that neighborhood there is no singularity. Dominated convergence justifies the limits and local continuity of both derivatives.

The function $u$ has a harmonic conjugate on the strip. With the normalization $v(\alpha,0)=0$, we can construct it by
\[
v(x,y):=-\int_\alpha^x u_y(s,0)\,ds
        +\int_0^y u_x(x,t)\,dt.
\]
The equality $u_{xx}+u_{yy}=0$ gives $v_y=u_x$ and $v_x=-u_y$. Thus $h:=u+iv$ is holomorphic, $h(\alpha)=u(\alpha,0)$ is real, and $\operatorname{Re}h(j+it)=\varphi_j(t)$.

Fix $\varepsilon,\delta>0$. The functions $t\mapsto\log(\|F(j+it)\|_{X_j}+\varepsilon)$ are continuous and tend to $\log\varepsilon$. Using a cutoff outside a large interval and convolution with a smooth function, choose $\varphi_j$ as above so that
\[
\log(\|F(j+it)\|_{X_j}+\varepsilon)
\leq\varphi_j(t)
\leq\log(\|F(j+it)\|_{X_j}+\varepsilon)+\delta.
\]
Construct $h$ with these data and define
\[
G(z):=e^{h(\alpha)-h(z)}F(z).
\]
Boundedness of $\operatorname{Re}h$ proves that $G\in\mathscr F$; moreover, $G(\alpha)=F(\alpha)$ and
\[
\|G(j+it)\|_{X_j}
\leq e^{u(\alpha,0)}
\frac{\|F(j+it)\|_{X_j}}{\|F(j+it)\|_{X_j}+\varepsilon}
\leq e^{u(\alpha,0)}.
\]
By the definition of the interpolated norm,
\[
\log\|F(\alpha)\|_{[X_0,X_1]_\alpha}
\leq\delta+\sum_{j=0}^1\int_{\mathbb R} P_j(\alpha,t)
 \log(\|F(j+it)\|_{X_j}+\varepsilon)\,dt.
\]
Let first $\delta\to0$ and then $\varepsilon\to 0^{+}$. Monotone convergence, applied after subtracting a constant upper bound, justifies the limit even when the integral equals $-\infty$.

To obtain \eqref{eq:poisson-potencias-norma-interpolada}, apply concavity of the logarithm with respect to the probability measures $m_j^{-1}P_j(\alpha,t)\,dt$ to the powers $r_j$. This integral inequality follows from $\log s\leq s-1$, first dividing by the average. If $r_j=\infty$, use the supremum. The same inequality with respect to the measure of total mass one on both lines gives \eqref{eq:poisson-lineal-norma-interpolada}.

Finally, for a function with bounded boundary values that do not
vanish at infinity, apply the preceding argument to
$F_d(z):=e^{d(z-\alpha)^2}F(z)$, with $d>0$.
This function belongs to $\mathscr F$ and $F_d(\alpha)=F(\alpha)$.
In the logarithmic estimate, the integral of the logarithm changes
by exactly
\[
d\sum_{j=0}^1\int_{\mathbb R}P_j(\alpha,t)
\bigl((j-\alpha)^2-t^2\bigr)\,dt.
\]
This quantity is finite and tends to zero with $d$, since the
Poisson kernels decay exponentially. The formula also holds if
the original integral equals $-\infty$: adding an integrable
function does not change that value. This gives the logarithmic
limit without requiring $\log\|F(j+it)\|_{X_j}$ to be integrable.
For finite powers and the linear bound, dominated convergence
does apply: if $0<d\leq1$, the boundary factor
$e^{d((j-\alpha)^2-t^2)}$ is bounded by $e$, and the original
norms are bounded. For an infinite power, use directly that this
factor is at most $e^d$.
\end{proof}

\subsection{Sums of spaces and Banach-valued functions}

Applications to Fourier blocks and families of charts require interpolation of several components at once. We first study a finite family, where each component admits a Calderón representative. We then pass to infinite families by truncation and to measurable functions by approximation with simple functions. The restrictions on the exponents express exactly when these approximations converge in norm.

For a Banach space $X$, write $L^r(A;X)$ for the classes of strongly measurable functions with finite norm
\[
\|f\|_{L^r(A;X)}:=\left(\int_A\|f(a)\|_X^r\,d\mu(a)\right)^{\frac{1}{r}}
\]
or with the essential supremum norm if $r=\infty$. These spaces are complete. Indeed, from a Cauchy sequence, extract a subsequence whose differences have summable norms. Minkowski's inequality and monotone convergence show that the pointwise sum of the norms of these differences belongs to $L^r(A)$ if $r<\infty$. Thus the series converges almost everywhere in $X$, and its tails tend to zero in $L^r(A;X)$. If $r=\infty$, the same conclusion holds outside the countable union of the null sets for the essential bounds. The limit is strongly measurable by Lemma~\ref{lem:b5-limite-fuertemente-medible}.

\begin{proposition}[Interpolation of $\ell^p(I)$-sums and Bochner spaces]
\label{prop:interpolacion-ellp-valores-banach}
Let $I$ be a set, let $1\leq r_0,r_1<\infty$ and $0<\theta<1$, and set
\[
\frac{1}{r}=\frac{1-\theta}{r_0}+\frac{\theta}{r_1}.
\]
For a family of compatible pairs of complex Banach spaces $(X_{0,i},X_{1,i})_{i\in I}$, we have isometrically
\begin{equation}
\left[\ell^{r_0}\bigl(I;(X_{0,i})_{i\in I}\bigr),
      \ell^{r_1}\bigl(I;(X_{1,i})_{i\in I}\bigr)\right]_\theta
=\ell^r\bigl(I;([X_{0,i},X_{1,i}]_\theta)_{i\in I}\bigr).
\label{eq:interpolacion-compleja-sumas-banach}
\end{equation}
If $w_{0,i},w_{1,i}>0$, the same identity with weighted norms is
\begin{equation}
\left[\ell^{r_0}\bigl(I;(w_{0,i}X_{0,i})_{i\in I}\bigr),
      \ell^{r_1}\bigl(I;(w_{1,i}X_{1,i})_{i\in I}\bigr)\right]_\theta
=\ell^r\bigl(I;(w_{0,i}^{1-\theta}w_{1,i}^\theta
      [X_{0,i},X_{1,i}]_\theta)_{i\in I}\bigr),
\label{eq:interpolacion-compleja-sumas-pesos}
\end{equation}
where $wX$ denotes $X$ with norm $w\|\cdot\|_X$. For a finite family, these identities also hold if some $r_j$ is infinite, with $\displaystyle\frac{1}{\infty}=0$.

If $(A,\mu)$ is $\sigma$-finite and $(X_0,X_1)$ is a compatible pair, then for $r_0,r_1<\infty$ we also have isometrically
\begin{equation}
[L^{r_0}(A;X_0),L^{r_1}(A;X_1)]_\theta
=L^r(A;[X_0,X_1]_\theta).
\label{eq:interpolacion-compleja-bochner}
\end{equation}
\end{proposition}

\begin{proof}
Begin with a finite family. Let $F$ be a Calderón function for the endpoint sums, and write $F_i$ for its $i$ coordinate. The coordinate projections are continuous at the endpoints, so $F_i\in\mathscr F(X_{0,i},X_{1,i})$. Apply \eqref{eq:poisson-potencias-norma-interpolada} to each $F_i$ and set
\[
b_{j,i}:=\left(\frac{1}{m_j}\int_{\mathbb R}P_j(\theta,t)
                  \|F_i(j+it)\|_{X_{j,i}}^{r_j}\,dt\right)^{\frac{1}{r_j}},
\quad m_0=1-\theta,\quad m_1=\theta,
\]
with suprema if $r_j=\infty$. We obtain
\[
\|F_i(\theta)\|_{[X_{0,i},X_{1,i}]_\theta}
\leq b_{0,i}^{1-\theta}b_{1,i}^\theta.
\]
If $r<\infty$, Hölder's inequality for sums, with $\displaystyle\frac{r(1-\theta)}{r_0}+r\displaystyle\frac{\theta}{r_1}=1$, gives
\[
\|(F_i(\theta))_i\|_{\ell^r([X_{0,i},X_{1,i}]_\theta)}
\leq\|(b_{0,i})_i\|_{\ell^{r_0}(I)}^{1-\theta}
      \|(b_{1,i})_i\|_{\ell^{r_1}(I)}^\theta
\leq\|F\|_{\mathscr F}.
\]
The last bound follows by interchanging sum and integral using Tonelli's theorem. If some $r_j$ is infinite, use the maximum of the $b_{j,i}$; if $r=\infty$, both endpoints are infinite and take the maximum in the coordinate estimate. This proves one inclusion.

For the other inclusion, let $x=(x_i)_i$ and set $a_i:=\|x_i\|_{[X_{0,i},X_{1,i}]_\theta}$ and $A:=\|(a_i)_i\|_{\ell^r(I)}$. If $A=0$, there is nothing to prove. Omit zero coordinates and choose representatives $F_i(\theta)=x_i$ with $\|F_i\|_{\mathscr F}\leq(1+\varepsilon)a_i$. When $r<\infty$, define
\[
\frac{1}{r(z)}:=\frac{1-z}{r_0}+\frac{z}{r_1},\qquad
F(z):=\left(A^{1-\frac{r}{r(z)}}a_i^{\frac{r}{r(z)}-1}F_i(z)\right)_i.
\]
This is a Calderón function for the sums, and $F(\theta)=x$. On each boundary,
\[
\|F(j+it)\|_{\ell^{r_j}(X_{j,i})}
\leq(1+\varepsilon)A^{1-\frac{r}{r_j}}
             \left(\sum_{i\in I} a_i^r\right)^{\frac{1}{r_j}}
=(1+\varepsilon)A,
\]
with the maximum interpretation if $r_j=\infty$. If both endpoints are infinite, simply take $F=(F_i)_i$; the family is finite. Letting $\varepsilon\to 0^{+}$ proves equality of norms.

For arbitrary $I$ and $r_0,r_1<\infty$, restriction to finite subsets gives the first inclusion and its bound, since the norm of $\ell^r(I)$ is the supremum of the norms of its finite restrictions. Conversely, every element of $\ell^r(I)$ has at most countably many nonzero coordinates: for each integer $k\geq1$, only finitely many coordinates have norm greater than $\displaystyle\frac{1}{k}$. Its truncations converge in $\ell^r(I)$, and the bound already proved makes them Cauchy in the interpolated space. Completeness and coordinate projections identify the limit. No density in an infinite $\ell^\infty(I)$-sum has been used.

The weighted formula is obtained first in each coordinate. Indeed, multiplying a representative by $w_0^{z-\theta}w_1^{\theta-z}$ preserves its value at $\theta$ and changes its weighted boundary norms by the common factor $w_0^{1-\theta}w_1^\theta$. Multiplication by the inverse gives the reverse inequality. Then apply the sum formula just proved.

Turn now to Bochner functions. For disjoint measurable sets $E_1,\ldots,E_N$ of positive finite measure, define
\[
P f:=\sum_{\nu=1}^N\mathbf1_{E_\nu}
                 \frac{1}{\mu(E_\nu)}\int_{E_\nu}f\,d\mu,
\]
and set $Pf=0$ outside their union. The Bochner integral inequality and Hölder's inequality show, for every Banach space $X$ and $1\leq s<\infty$,
\[
\|Pf\|_{L^s(A;X)}\leq\|f\|_{L^s(A;X)}.
\]
Its image is identified with a finite sum of copies of $X$ whose norm has weights $\mu(E_\nu)^{\displaystyle\frac{1}{s}}$. The weighted formula gives exactly the weight $\mu(E_\nu)^{\displaystyle\frac{1}{r}}$ in the interpolated space.

Let $F\in\mathscr F(L^{r_0}(A;X_0),L^{r_1}(A;X_1))$. For each $j$, the image of the boundary line together with zero is compact in $L^{r_j}(A;X_j)$. Choose a finite net in that image and approximate its elements by simple functions supported on sets of finite measure. Such approximations exist by strong measurability and integrability of the norm: first truncate the domain and the size of the function, then use the simple functions in its definition. A common partition of the supports of these approximants produces an operator $P$ fixing them. Its contractivity implies, for each function $f$ in the net and its approximation $g$,
\[
\|Pf-f\|\leq\|P(f-g)\|+\|g-f\|\leq2\|f-g\|.
\]
Taking sufficiently fine nets and small errors gives operators $P_k$ such that $P_kF\to F$ in $\mathscr F$.

The equality for finite sums implies
\[
\|P_kF(\theta)-P_lF(\theta)\|_{L^r(A;[X_0,X_1]_\theta)}
\leq\|(P_k-P_l)F\|_{\mathscr F}.
\]
To use it, simply refine both partitions to a common partition. Thus $P_kF(\theta)$ converges in the Bochner space on the right-hand side of \eqref{eq:interpolacion-compleja-bochner}. That limit agrees with $F(\theta)$: both convergences imply local convergence in measure with values in $X_0+X_1$, by Hölder's inequality on sets of finite measure. This proves the inclusion from the left-hand side. The bound is one, since $\|P_kF\|_{\mathscr F}\leq\|F\|_{\mathscr F}$.

Finally, simple functions supported on sets of finite measure with values in $[X_0,X_1]_\theta$ are dense in their Bochner space. For each such function, the construction using coordinate representatives and the preceding weights gives the reverse inclusion with bound one. Differences of two simple functions are treated on a common partition; completeness allows passage to the limit. This completes the proof.
\end{proof}

\subsection{Complex reiteration}

To prove reiteration while retaining the general assumptions, we shall use a description through Fourier coefficients of periodic functions. The idea is to express an element in two ways: as the sum of a sequence and as the constant value of a sequence decomposed between the endpoints. The relationship between them is obtained by taking differences and partial sums. This procedure is related to the approach of \cite[\S\S4 and 6]{LindemulderLorist2024}.

\begin{definition}[Fourier coefficient spaces]
\label{def:coeficientes-fourier-interpolacion}
Let $X$ be a complex Banach space. Denote by $\mathcal C(X)$ the space of sequences $u=(u_n)_{n\in\mathbb Z}$ of the form
\[
u_n=\frac{1}{2\pi}\int_{-\pi}^\pi f(t)e^{-int}\,dt,
\qquad f\in L^1((-\pi,\pi);X),
\]
with the endpoints of the interval identified. Its norm is
\[
\|u\|_{\mathcal C(X)}:=\frac{1}{2\pi}\int_{-\pi}^\pi\|f(t)\|_X\,dt.
\]
For $b>1$ and $\rho\in\mathbb R$, write
\[
\mathcal C_b(\rho;X):=\{u\mid(b^{\rho n}u_n)_n\in\mathcal C(X)\},
\qquad
\|u\|_{\mathcal C_b(\rho;X)}:=\|(b^{\rho n}u_n)_n\|_{\mathcal C(X)}.
\]
These periodic coefficients use the measure $\displaystyle\frac{dt}{2\pi}$; they are not the Fourier transform on $\mathbb R$ defined in the preceding chapter.
\end{definition}

We justify that the norm is well defined. The Fejér kernel
\[
\Phi_N(t):=\frac{1}{N+1}\left|\sum_{k=0}^N e^{ikt}\right|^2
\]
is nonnegative, has integral one with respect to $\displaystyle\frac{dt}{2\pi}$, and satisfies $\Phi_N(t)\leq((N+1)\sin^2(\displaystyle\frac{t}{2}))^{-1}$ outside $t=0$. Its convolutions therefore approximate continuous periodic functions uniformly: use uniform continuity near zero and the preceding bound on the complement. Approximating a Bochner function in $L^1((-\pi,\pi);X)$ by simple functions and their scalar indicator functions by continuous periodic functions gives the same result in $L^1((-\pi,\pi);X)$, since the convolution operators have norm at most one. Their coefficients are
\[
(Q_Nu)_n:=\left(1-\frac{|n|}{N+1}\right)_+u_n.
\]
Consequently, if all coefficients of $f$ vanish, then $\Phi_N*f=0$ for every $N$ and $f=0$. This proves uniqueness of $f$ and completeness of $\mathcal C(X)$ through its isometric identification with $L^1((-\pi,\pi);X)$ with normalized measure. It also proves
\begin{equation}
\|Q_Nu\|_{\mathcal C_b(\rho;X)}\leq\|u\|_{\mathcal C_b(\rho;X)},
\qquad Q_Nu\longrightarrow u\quad\text{in }\mathcal C_b(\rho;X).
\label{eq:fejer-coeficientes-interpolacion}
\end{equation}
Each coordinate satisfies $\|b^{\rho n}u_n\|_X\leq
\|u\|_{\mathcal C_b(\rho;X)}$. Thus the coefficient spaces form compatible pairs in the space of sequences with values in $X_0+X_1$, equipped with the topology of coordinatewise convergence.

\begin{lemma}[Discrete representations of the complex method]
\label{lem:representaciones-discretas-complejas}
Let $b>1$ and $0<\theta<1$. For a compatible pair $(X_0,X_1)$, set
\[
\mathcal U_0:=\mathcal C_b(-\theta;X_0),\qquad
\mathcal U_1:=\mathcal C_b(1-\theta;X_1).
\]
Then $x\in[X_0,X_1]_\theta$ if and only if there is $u\in\mathcal U_0\cap\mathcal U_1$ such that
\[
x=\sum_{n\in\mathbb Z}u_n\quad\text{in }X_0+X_1.
\]
Its norm is equivalent to the infimum of $\displaystyle \max\{\|u\|_{\mathcal U_0},\|u\|_{\mathcal U_1}\}$ over these representations. It is also equivalent to
\[
\|\underline x\|_{\mathcal U_0+\mathcal U_1},
\qquad \underline x:=(x)_{n\in\mathbb Z}.
\]
In particular, the last quantity is finite exactly for elements of $[X_0,X_1]_\theta$. The constants depend only on $b$ and $\theta$.
\end{lemma}

\begin{proof}
If $u\in\mathcal U_0\cap\mathcal U_1$, the coordinate bounds give
\[
\|u_n\|_{X_0}\leq b^{\theta n}\|u\|_{\mathcal U_0},\qquad
\|u_n\|_{X_1}\leq b^{-(1-\theta)n}\|u\|_{\mathcal U_1}.
\]
Use the first for $n<0$ and the second for $n\geq0$. The geometric series prove absolute convergence in $X_0+X_1$. The exponential polynomials
\[
H_N(z):=\sum_{n\in \mathbb Z} b^{(z-\theta)n}(Q_Nu)_n
\]
take values in $X_0\cap X_1$ and have imaginary period $P:=\displaystyle\frac{2\pi}{\log b}$. On the boundary $j+it$, the corresponding periodic function is the Fejér sum of the function with coefficients $b^{(j-\theta)n}u_n$.

Lemma~\ref{lem:poisson-norma-interpolada}, in its version for bounded functions, implies
\begin{equation}
\|H_N(\theta)-H_M(\theta)\|_{[X_0,X_1]_\theta}
\leq C_{b,\theta}\max_{j\in\{0,1\}}\|(Q_N-Q_M)u\|_{\mathcal U_j}.
\label{eq:estimacion-polinomios-periodicos}
\end{equation}
To check the constant, split each Poisson integral into intervals of length $P$. The periodized kernel $\displaystyle \sum_{k\in\mathbb Z}P_j(\theta,t+kP)$ is uniformly bounded on $[0,P]$, since $P_j$ is continuous and decays exponentially. This bounds the integral by a constant times the average $\displaystyle P^{-1}\int_0^P\|H_N(j+it)-H_M(j+it)\|_{X_j}\,dt$. By \eqref{eq:fejer-coeficientes-interpolacion}, the values $H_N(\theta)$ converge in the interpolated space; in the sum space they converge to $\displaystyle \sum_{n\in \mathbb Z} u_n$. This gives the first norm bound.

We prove the converse direction. If $H$ is holomorphic on the strip, continuous and bounded up to the boundary in the corresponding norms, and has imaginary period $P$, define
\begin{equation}
u_n:=b^{-(j-\theta)n}\frac{1}{P}\int_0^P
 H(j+it)e^{-in(\log b)t}\,dt,
\qquad j=0,1.
\label{eq:coeficientes-representante-periodico}
\end{equation}
The two expressions agree: apply Cauchy's theorem to $b^{-n(z-\theta)}H(z)$ on a rectangle of height $P$; the horizontal integrals cancel by periodicity. Continuity allows passage to the boundary lines. Thus $u_n\in X_0\cap X_1$, $u\in\mathcal U_0\cap\mathcal U_1$, and
\[
\max_{j\in\{0,1\}}\|u\|_{\mathcal U_j}
\leq\max_{j\in\{0,1\}}\sup_{t\in\mathbb R}\|H(j+it)\|_{X_j}.
\]
The same integrals on the line $\operatorname{Re}z=\theta$ show that $(u_n)_n$ are the coefficients of $H(\theta+it)$. Fejér convergence on that line and the absolute convergence already proved give $\displaystyle H(\theta)=\displaystyle\sum_{n\in\mathbb Z}u_n$.

Now let $x\in[X_0,X_1]_\theta$, $x\neq0$, and choose $F(\theta)=x$ with $\|F\|_{\mathscr F}\leq2\|x\|_{[X_0,X_1]_\theta}$. Fix $d>0$ sufficiently large that
\[
2\sum_{k\in\mathbb Z\setminus\{0\}}e^{-dk^2P^2}\leq\frac{1}{2}
\]
and periodize:
\[
H(z):=\sum_{k\in\mathbb Z}
 e^{d(z+ikP-\theta)^2}F(z+ikP).
\]
The series converges uniformly on compact subsets and on one period of each boundary, with values in $X_j$ on the boundary $j$. Moreover,
\[
\max_{j\in\{0,1\}}\sup_{t\in\mathbb R}\|H(j+it)\|_{X_j}
\leq C_{b,\theta}\|F\|_{\mathscr F},
\]
since $\displaystyle \sup_{t\in\mathbb R}\displaystyle\sum_{k\in \mathbb Z} e^{-d(t+kP)^2}<\infty$. If $y:=H(\theta)$, Lemma~\ref{lem:regularizacion-gaussiana-calderon} and the choice of $d$ give
\[
\|x-y\|_{[X_0,X_1]_\theta}
\leq\sum_{k\in\mathbb Z\setminus\{0\}}e^{-dk^2P^2}
       \|F(\theta+ikP)\|_{[X_0,X_1]_\theta}
\leq\frac{1}{2}\|x\|_{[X_0,X_1]_\theta}.
\]
Formula \eqref{eq:coeficientes-representante-periodico} represents $y$ with norm in $\mathcal U_0\cap\mathcal U_1$ at most $C_{b,\theta}\|x\|_{[X_0,X_1]_\theta}$. Repeat the construction for $x-y$ and successive remainders. Their norms are at most $2^{-l}\|x\|$. The resulting coefficient sequences form an absolutely convergent series in the Banach space $\mathcal U_0\cap\mathcal U_1$. Its sum $u$ satisfies $\displaystyle \sum_{n\in\mathbb Z}u_n=x$ and has the required bound. Interchange of the two sums is justified by continuity of the operator $\displaystyle u\mapsto\displaystyle\sum_{n\in\mathbb Z}u_n$ into the sum space, proved using geometric series at the start.

It remains to establish the characterization through $\underline x$. Given a representation $\displaystyle x=\displaystyle\sum_{n\in\mathbb Z}u_n$, set
\[
a_n:=\sum_{\substack{k\in\mathbb Z\\k<n}}u_k\in X_0,\qquad
c_n:=\sum_{\substack{k\in\mathbb Z\\k\geq n}}u_k\in X_1.
\]
Then $a_n+c_n=x$. For the weighted coefficients,
\[
\begin{aligned}
b^{-\theta n}a_n
 &=\sum_{l\in\mathbb N}b^{-\theta l}
                  b^{-\theta(n-l)}u_{n-l},\\
b^{(1-\theta)n}c_n
 &=\sum_{l\in\mathbb N_0}b^{-(1-\theta)l}
                  b^{(1-\theta)(n+l)}u_{n+l}.
\end{aligned}
\]
Shifting the indices of a sequence in $\mathcal C(X)$ corresponds to multiplying its function by $e^{ilt}$, an isometry on $L^1((-\pi,\pi);X)$. The preceding geometric series give
\[
\|a\|_{\mathcal U_0}+\|c\|_{\mathcal U_1}
\leq C_{b,\theta}\max_{j\in\{0,1\}}\|u\|_{\mathcal U_j}.
\]
Conversely, if $x=a_n+c_n$ with $a\in\mathcal U_0$ and $c\in\mathcal U_1$, define
\[
u_n:=a_{n+1}-a_n=c_n-c_{n+1}\in X_0\cap X_1.
\]
The same shift invariance implies
\[
\|u\|_{\mathcal U_0}\leq(1+b^\theta)\|a\|_{\mathcal U_0},\qquad
\|u\|_{\mathcal U_1}\leq(1+b^{-(1-\theta)})\|c\|_{\mathcal U_1}.
\]
The coordinate bounds show that $a_n\to0$ in $X_0$ as $n\to-\infty$, and $c_n\to0$ in $X_1$ as $n\to\infty$. Thus $\displaystyle \sum_{n=-N}^Nu_n=a_{N+1}-a_{-N}\to x$ in the sum space. Taking infima completes the proof.
\end{proof}

\begin{lemma}[Decomposition of coefficients at an intermediate level]
\label{lem:coeficientes-nivel-intermedio}
Let $b>1$ and $0\leq\alpha\leq1$, and set $E_\alpha:=[X_0,X_1]_\alpha$ if $0<\alpha<1$, $E_0:=X_0$, and $E_1:=X_1$. Then
\begin{equation}
\mathcal C(E_\alpha)
\hookrightarrow
\mathcal C_b(-\alpha;X_0)+\mathcal C_b(1-\alpha;X_1)
\label{eq:coeficientes-encaje-suma}
\end{equation}
with norm at most one.
\end{lemma}

\begin{proof}
If $\alpha=0$ or $1$, the space on the left equals one of the summands, and the assertion is immediate. Assume $0<\alpha<1$ and denote the sum space on the right by $V$. Begin with a trigonometric polynomial
\[
f(t)=\sum_{\substack{n\in\mathbb Z\\|n|\leq N}}u_ne^{int},\qquad u_n\in X_0\cap X_1.
\]
Thus its coefficient sequence belongs to $V$. Let $\Lambda\in V'$ have norm at most one. By \eqref{eq:fejer-coeficientes-interpolacion}, the functionals $\Lambda\circ Q_M$ have norm at most one and converge to $\Lambda$ on each element of $V$. Each depends on only finitely many coordinates, so it has the form
\[
(\Lambda\circ Q_M)(v)=\sum_{n\in \mathbb Z} y_n(v_n),
\qquad y_n\in(X_0+X_1)',
\]
with only finitely many nonzero $y_n$. For $j\in\{0,1\}$, set
\[
H(z,t):=\sum_{n\in \mathbb Z} b^{(\alpha-z)n}e^{-int}y_n.
\]
The norm of $\Lambda\circ Q_M$ restricted to $\mathcal C_b(j-\alpha;X_j)$ is
\begin{equation}
\sup_{t\in[-\pi,\pi]}\|H(j,t)\|_{X_j'}\leq1.
\label{eq:norma-funcional-coeficientes}
\end{equation}
To verify the equality giving this norm, if $(b^{(j-\alpha)n}v_n)_n$ are the coefficients of $g\in L^1((-\pi,\pi);X_j)$, then
\[
(\Lambda\circ Q_M)(v)=\frac{1}{2\pi}\int_{-\pi}^\pi H(j,t)(g(t))\,dt.
\]
This proves one inequality. For the other, if $X_j=\{0\}$,
both norms are zero. Otherwise, choose $t_0$ and a vector in the
unit ball of $X_j$ whose value under $H(j,t_0)$ approximates its
dual norm, and multiply it by a phase to make that value real and
nonnegative. Multiply the vector by the indicator of a small
interval around $t_0$, normalized in $L^1((-\pi,\pi))$.
Continuity of $H(j,\cdot)$ allows the interval length to tend to
zero and recovers $\|H(j,t_0)\|_{X_j'}$. Taking the supremum
over $t_0$ proves equality of norms.

If $F\in\mathscr F(X_0,X_1)$ and $t$ is fixed, the scalar function $z\mapsto H(z,t)(F(z))$ is holomorphic and bounded on the strip. On the boundary $j+is$,
\[
H(j+is,t)=H(j,t+s\log b),
\]
so \eqref{eq:norma-funcional-coeficientes} and the three-lines lemma imply
\[
|H(\alpha,t)(F(\alpha))|\leq\|F\|_{\mathscr F}.
\]
Taking the infimum over $F$ shows that $H(\alpha,t)$ has norm at most one on $E_\alpha$. Consequently,
\[
|(\Lambda\circ Q_M)(u)|
=\left|\frac{1}{2\pi}\int_{-\pi}^\pi H(\alpha,t)(f(t))\,dt\right|
\leq\frac{1}{2\pi}\int_{-\pi}^\pi\|f(t)\|_{E_\alpha}\,dt.
\]
Let $M\to\infty$ and take the supremum over $\Lambda$. Hahn--Banach gives $\|u\|_V\leq\|u\|_{\mathcal C(E_\alpha)}$.

Lemma~\ref{lem:regularizacion-gaussiana-calderon} ensures that $X_0\cap X_1$ is dense in $E_\alpha$. Hence trigonometric polynomials with coefficients in that intersection are dense in $L^1((-\pi,\pi);E_\alpha)$: first apply Fejér approximation, then approximate its finitely many coefficients. The bound obtained and completeness of $V$ extend the embedding to all of $\mathcal C(E_\alpha)$. Coordinate projections show that this extension preserves each coefficient, so it is the inclusion in \eqref{eq:coeficientes-encaje-suma}.
\end{proof}

\begin{theorem}[Complex reiteration]
\label{teo:reiteracion-interpolacion-compleja}
Let $0\leq\theta_0<\theta_1\leq1$, $0<\eta<1$, and $\theta:=(1-\eta)\theta_0+\eta\theta_1$. If $X_{\theta_j}:=[X_0,X_1]_{\theta_j}$ at interior points, with $X_0,X_1$ retained at the endpoints, then
\[
[X_{\theta_0},X_{\theta_1}]_\eta=[X_0,X_1]_\theta
\]
with equivalent norms. The constants depend only on $\theta_0,\theta_1,\eta$.
\end{theorem}

\begin{proof}
Set $Y_j:=X_{\theta_j}$ and $\delta:=\theta_1-\theta_0$. For the first inclusion, let $F\in\mathscr F(X_0,X_1)$ and consider its regularizations $F_\varepsilon$ from Lemma~\ref{lem:regularizacion-gaussiana-calderon}. The functions
\[
G_\varepsilon(z):=F_\varepsilon(\theta_0+\delta z)
\]
are holomorphic on $X_0\cap X_1$ and belong to $\mathscr F(Y_0,Y_1)$. The bound on interior lines in the same lemma gives
\[
\|G_\varepsilon\|_{\mathscr F(Y_0,Y_1)}
\leq\|F_\varepsilon\|_{\mathscr F(X_0,X_1)}
\leq\|F\|_{\mathscr F(X_0,X_1)}.
\]
On lines coinciding with an endpoint, use the boundary norm directly. To check boundedness in the sum space and the required holomorphicity, use the fact that $F_\varepsilon$ is entire in the intersection and formula \eqref{eq:regularizacion-calderon-dos-fronteras}: its norms at both endpoints are bounded on every vertical strip of finite width. The same estimate applied to $F_\varepsilon-F_{\varepsilon'}$ shows that $(G_\varepsilon)$ is Cauchy in $\mathscr F(Y_0,Y_1)$. Its limit takes the value $F(\theta)$ at $\eta$ by continuity of the inclusions into $X_0+X_1$. Thus,
\[
\|F(\theta)\|_{[Y_0,Y_1]_\eta}\leq\|F\|_{\mathscr F(X_0,X_1)}.
\]

No density assumption at a retained endpoint is needed: in
Lemma~\ref{lem:coeficientes-nivel-intermedio}, the case $\alpha=0$
is the inclusion $w\mapsto w+0$ into the sum of coefficient spaces,
and the case $\alpha=1$ is $w\mapsto0+w$.

For the reverse inclusion, fix $b=e$ and $c:=b^\delta>1$. Let $x\in[Y_0,Y_1]_\eta$. The constant-sequence characterization in Lemma~\ref{lem:representaciones-discretas-complejas}, applied to the pair $(Y_0,Y_1)$ with base $c$, allows us to write
\[
\underline x=v_0+v_1,\qquad
v_j\in\mathcal C_c(j-\eta;Y_j),
\]
with sum of norms at most $C\|x\|_{[Y_0,Y_1]_\eta}$, up to arbitrarily small error. Since $\delta(j-\eta)=\theta_j-\theta$, the sequences
\[
w_{j,n}:=b^{(\theta_j-\theta)n}v_{j,n}
\]
belong to $\mathcal C(Y_j)$ with the same norm as $v_j$ in its weighted space. Lemma~\ref{lem:coeficientes-nivel-intermedio} decomposes them as
\[
w_j=a_{j,0}+a_{j,1},\qquad
a_{j,0}\in\mathcal C_b(-\theta_j;X_0),\qquad
a_{j,1}\in\mathcal C_b(1-\theta_j;X_1),
\]
with sum of norms at most $\|w_j\|_{\mathcal C(Y_j)}$, up to arbitrarily small error. Multiplying each coordinate by $b^{(\theta-\theta_j)n}$ and summing over $j\in\{0,1\}$ gives
\[
\underline x=a_0+a_1,\qquad
a_0\in\mathcal C_b(-\theta;X_0),\qquad
a_1\in\mathcal C_b(1-\theta;X_1).
\]
The weights cancel exactly: for example, $b^{-\theta n}b^{(\theta-\theta_j)n}=b^{-\theta_jn}$. Therefore,
\[
\|\underline x\|_{\mathcal C_b(-\theta;X_0)+\mathcal C_b(1-\theta;X_1)}
\leq C\|x\|_{[Y_0,Y_1]_\eta}.
\]
Another application of Lemma~\ref{lem:representaciones-discretas-complejas} implies $x\in[X_0,X_1]_\theta$ and gives the reverse bound. The constants depend only on the bases $b,c$ and the indicated parameters; the decomposition lemma has norm one. This includes $\theta_0=0$ or $\theta_1=1$.
\end{proof}

\begin{lemma}[Calderón's gluing lemma]
\label{lem:pegado-calderon}
Let $0\leq\alpha_0<\alpha_1\leq1$, let $0<\eta<1$, and set $\alpha:=(1-\eta)\alpha_0+\eta\alpha_1$. If $G\in\mathscr F(X_{\alpha_0},X_{\alpha_1})$, then for each $\varepsilon>0$ there is $F_\varepsilon\in\mathscr F(X_0,X_1)$ such that
\[
F_\varepsilon(\alpha)=G(\eta),\qquad
\|F_\varepsilon\|_{\mathscr F(X_0,X_1)}
\leq C_{\alpha_0,\alpha_1,\eta}
             \bigl(\|G\|_{\mathscr F(X_{\alpha_0},X_{\alpha_1})}
                   +\varepsilon\bigr).
\]
Also, for $F\in\mathscr F(X_0,X_1)$ and $0<\alpha<1$,
\[
\sup_{t\in\mathbb R}\|F(\alpha+it)\|_{[X_0,X_1]_\alpha}\leq\|F\|_{\mathscr F}.
\]
\end{lemma}

\begin{proof}
Theorem~\ref{teo:reiteracion-interpolacion-compleja} gives $G(\eta)\in[X_0,X_1]_\alpha$ with norm at most $C\|G\|_{\mathscr F}$. Choose a representative whose value at $\alpha$ is $G(\eta)$ and whose norm differs by less than $\varepsilon$ from the infimum defining that norm. The final assertion is Lemma~\ref{lem:regularizacion-gaussiana-calderon}.
\end{proof}

\begin{proposition}[Real interpolation of $\ell^p(I)$-sums]
\label{prop:interpolacion-real-ellp}
Let $1\leq p<\infty$. For every compatible pair $(X_0,X_1)$ and $0<\theta<1$,
\[
(\ell^p(X_0),\ell^p(X_1))_{\theta,p}
=
\ell^p((X_0,X_1)_{\theta,p})
\]
as sets, and for every sequence $x=(x_j)$,
\[
 \|x\|_{\ell^p((X_0,X_1)_{\theta,p})}
 \leq
 \|x\|_{(\ell^p(X_0),\ell^p(X_1))_{\theta,p}}
 \leq
 2^{1-\frac{1}{p}}
 \|x\|_{\ell^p((X_0,X_1)_{\theta,p})}.
\]
\end{proposition}

\begin{proof}
For a sequence $x=(x_j)_{j\in\mathbb N}$, set $J:=\mathbb N$ and
\[
 A(t,x):=\left(\sum_{j\in J}
 K(t,x_j;X_0,X_1)^p\right)^{\frac{1}{p}}.
\]
If $x=a+b$, Minkowski's inequality from Theorem~\ref{teo:b6-espacios-lp-espacios-de-lebesgue}, applied to the discrete measure space $J$, gives
\[
 A(t,x)\leq\|a\|_{\ell^p(X_0)}+t\|b\|_{\ell^p(X_1)}.
\]
Taking the infimum yields $A(t,x)\leq K(t,x;\ell^p(X_0),\ell^p(X_1))$. Conversely, for each coordinate choose $x_j=a_j+b_j$ with $\|a_j\|_{X_0}+t\|b_j\|_{X_1}\leq K(t,x_j)+\varepsilon_j$, where $(\varepsilon_j)\in\ell^p(J)$ is arbitrarily small. Then
\[
 \|a\|_{\ell^p(X_0)}+t\|b\|_{\ell^p(X_1)}
 \leq2^{1-\frac{1}{p}}
 \left(\sum_{j\in J}
  (\|a_j\|_{X_0}+t\|b_j\|_{X_1})^p\right)^{\frac{1}{p}}.
\]
Letting $(\varepsilon_j)$ tend to zero gives
\[
 A(t,x)\leq K(t,x;\ell^p(X_0),\ell^p(X_1))
 \leq2^{1-\frac{1}{p}}A(t,x).
\]
Raise to the power $p$, integrate with respect to $t^{-\theta p}\displaystyle\frac{dt}{t}$, and apply Tonelli's theorem, Theorem~\ref{teo:tonelli}. The estimates also hold if either side is infinite. If the right-hand side is finite, the estimate at $t=1$ and the chosen decompositions show that $x$ belongs to the sum space; hence both inclusions have been justified.
\end{proof}

\section{Lorentz spaces as interpolation spaces}
\label{sec:lorentz-interpolacion}

Fractional embeddings obtained by the real method naturally take values in Lorentz spaces. We recall only what will be used later.

\begin{definition}[Decreasing rearrangement and Lorentz spaces]
\label{def:espacios-lorentz}\index{Lorentz space}
Let $(A,\mu)$ be a $\sigma$-finite nonatomic measure space, and let $f\colon A\longrightarrow\mathbb C$ be measurable. Define
\[
\mu_f(\lambda):=\mu(\{x\in A\mid |f(x)|>\lambda\}),\qquad
f^*(t):=\inf\{\lambda\in(0,\infty)\mid\mu_f(\lambda)\leq t\},\qquad t>0,
\]
with $\displaystyle \inf\varnothing=\infty$. For $1<p<\infty$ and $1\leq q\leq\infty$, the space $L^{p,q}(A)$ consists of the classes of functions for which
\[
\|f\|_{L^{p,q}(A)}
:=\|t^{\frac{1}{p}}f^*(t)\|_{L^q((0,\infty),\frac{dt}{t})}<\infty,
\]
and, if \(q=\infty\), define explicitly \[
\|f\|_{L^{p,\infty}(A)}
=\sup_{t\in(0,\infty)}t^{\frac{1}{p}}f^*(t).
\] We shall also use
\[
f^{**}(t):=\frac{1}{t}\int_0^t f^*(s)\,ds,\qquad
\|f\|_{L^{p,q}(A)}^\#:=\|t^{\frac{1}{p}}f^{**}(t)\|_{L^q((0,\infty),\frac{dt}{t})}.
\]
The first expression is a quasinorm; the second will be an equivalent norm, as we shall prove below.
\end{definition}

\begin{lemma}[Rearrangement identities]
\label{lem:identidades-reordenamiento-decreciente}
Under the assumptions of the preceding definition, for $t>0$ we have
\[
\int_0^t f^*(s)\,ds
=\sup_{\substack{E\subseteq A\text{ measurable}\\\mu(E)\leq t}}
  \int_E|f|\,d\mu.
\]
Consequently, if $f=g+h$,
\[
\int_0^t f^*(s)\,ds
\leq\int_0^t g^*(s)\,ds+\int_0^t h^*(s)\,ds.
\]
Moreover, for $0<p<\infty$,
\[
\int_A|f|^p\,d\mu=\int_0^\infty f^*(t)^p\,dt.
\]
\end{lemma}

\begin{proof}
The definition of the generalized inverse and right continuity of $\mu_f$ give
\[
|\{s>0\mid f^*(s)>\lambda\}|=\mu_f(\lambda),\qquad\lambda>0.
\]
Indeed, $s\geq\mu_f(\lambda)$ implies $f^*(s)\leq\lambda$. If $s<\mu_f(\lambda)$, right continuity provides $\delta>0$ with $\mu_f(\lambda+\delta)>s$, and then $f^*(s)\geq\lambda+\delta$. This argument includes the case $\mu_f(\lambda)=\infty$, by continuity from below of measures as the level decreases to $\lambda$. Cavalieri's formula from Proposition~\ref{prop:cavalieri-chebyshev} therefore gives
\begin{equation}
\int_0^t f^*(s)\,ds
=\int_0^\infty\min\{t,\mu_f(\lambda)\}\,d\lambda.
\label{eq:cavalieri-reordenamiento-truncado}
\end{equation}
For $\mu(E)\leq t$, Cavalieri's formula also implies
\[
\int_E|f|\,d\mu
=\int_0^\infty\mu(E\cap\{|f|>\lambda\})\,d\lambda
\leq\int_0^\infty\min\{t,\mu_f(\lambda)\}\,d\lambda.
\]

For the reverse inequality, we shall use the fact that in a nonatomic space, a measurable set of finite measure contains subsets of every measure between zero and its own measure. Here is a justification. Every set of positive measure contains a subset of positive measure less than any $\varepsilon>0$: split it into two parts of positive measure and repeat with the smaller part. To obtain a prescribed measure $a$, choose disjoint parts without exceeding $a$ and, at each step, with measure greater than half the supremum of still admissible measures. If the sum of their measures tended to a value less than $a$, the remainder would contain a fixed part of positive measure smaller than the difference; it would remain admissible at every step, contradicting convergence of the sum. The union of the parts has measure $a$.

If $|f|$ is simple and supported on a set of finite measure, choose $E$ by taking its highest level sets first and, when necessary, a portion of the last level to reach measure $t$. This attains the quantity in \eqref{eq:cavalieri-reordenamiento-truncado}. For a general function, approximate $|f|$ by an increasing sequence of nonnegative simple functions supported on sets of finite measure. Their distribution functions increase to $\mu_f$; monotone convergence in \eqref{eq:cavalieri-reordenamiento-truncado} allows passage to the limit and proves equality of the suprema.

Applying that equality to $|g+h|\leq|g|+|h|$ gives subadditivity. Finally, $|f|$ and $f^*$ have the same distribution function; Cavalieri's formula applied to their $p$th powers gives the final identity.
\end{proof}

The preceding subadditivity and Minkowski's inequality show that $\|\cdot\|_{L^{p,q}(A)}^\#$ is a norm. Since $f^*\leq f^{**}$, one comparison is immediate. For the other, set $g(t):=t^{\displaystyle\frac{1}{p}}f^*(t)$ and observe that
\[
t^{\frac{1}{p}}f^{**}(t)
=\int_0^t\left(\frac{s}{t}\right)^{1-\frac{1}{p}}g(s)\,\frac{ds}{s}.
\]
Theorem~\ref{teo:desigualdades-hardy-integrales}, with $1-\displaystyle\frac{1}{p}=\displaystyle\frac{1}{p'}>0$, gives
\begin{equation}
\|f\|_{L^{p,q}(A)}\leq\|f\|_{L^{p,q}(A)}^\#\leq p'\|f\|_{L^{p,q}(A)}.
\label{eq:normas-lorentz-equivalentes}
\end{equation}
This also proves the quasi-triangle inequality for the first expression. The next theorem will provide completeness for the equivalent norm.

\begin{theorem}[Interpolation of $L^p(A)$ spaces]\label{teo:interpolacion-Lp-lorentz}
Let $(A,\mu)$ be a $\sigma$-finite nonatomic measure space. Let $1\leq p_0<p_1\leq\infty$, $0<\theta<1$, $1\leq q\leq\infty$, and
\[
\frac{1}{p}=\frac{1-\theta}{p_0}+\frac{\theta}{p_1}.
\]
Then
\[
(L^{p_0}(A),L^{p_1}(A))_{\theta,q}
=
L^{p,q}(A)
\]
as sets. More precisely, there are $c_{p_0,p_1,\theta,q},C_{p_0,p_1,\theta,q}>0$, independent of $f$ and $(A,\mu)$, such that
\[
 c_{p_0,p_1,\theta,q}\|f\|_{L^{p,q}(A)}
 \leq\|f\|_{(L^{p_0}(A),L^{p_1}(A))_{\theta,q}}
 \leq C_{p_0,p_1,\theta,q}\|f\|_{L^{p,q}(A)}.
\]
Moreover, $L^{p,p}(A)=L^p(A)$ with equality of norms, and if $q_0\leq q_1$, then $L^{p,q_0}(A)\hookrightarrow L^{p,q_1}(A)$.
\end{theorem}

\begin{proof}
We first prove the identity
\begin{equation}
K(t,f;L^1(A),L^\infty(A))=\int_0^t f^*(s)\,ds.
\label{eq:funcional-K-L1-Linfty}
\end{equation}
If $f=g+h$, Lemma~\ref{lem:identidades-reordenamiento-decreciente} gives
\[
\int_0^t f^*(s)\,ds\leq\|g\|_{L^1(A)}+t\|h\|_{L^\infty(A)}.
\]
Taking the infimum proves one inequality. For the other, suppose that the integral is finite and set $\lambda:=f^*(t)<\infty$. With $\operatorname{sgn}f=\displaystyle\frac{f}{|f|}$ where $f\neq0$, and zero at its zeros, define
\[
g:=\operatorname{sgn}f\,(|f|-\lambda)_+,\qquad h:=f-g.
\]
Monotonicity of $\mu_f$ and the definition of $\lambda$ imply $\mu_f(s)>t$ if $0<s<\lambda$ and $\mu_f(s)\leq t$ if $s>\lambda$. By Cavalieri's formula and \eqref{eq:cavalieri-reordenamiento-truncado},
\[
\|g\|_{L^1(A)}+t\|h\|_{L^\infty(A)}
\leq\int_\lambda^\infty\mu_f(s)\,ds+t\lambda
=\int_0^t f^*(s)\,ds.
\]
In particular, this construction proves $f\in L^1(A)+L^\infty(A)$ whenever the integral is finite. If it is infinite, the first inequality already precludes a decomposition of finite cost. This gives \eqref{eq:funcional-K-L1-Linfty} with the extended interpretation.

For $1<p<\infty$, set $\alpha:=1-\displaystyle\frac{1}{p}$. The preceding identity gives
\[
\|f\|_{(L^1(A),L^\infty(A))_{\alpha,q}}
=\|t^{\frac{1}{p}}f^{**}(t)\|_{L^q((0,\infty),\frac{dt}{t})}=\|f\|_{L^{p,q}(A)}^\#.
\]
The comparisons \eqref{eq:normas-lorentz-equivalentes} prove $(L^1(A),L^\infty(A))_{1-\displaystyle\frac{1}{p},q}=L^{p,q}(A)$, as well as completeness with the norm $\|\cdot\|^\#$, by Theorem~\ref{teo:propiedades-basicas-metodo-K}.

For general exponents, write $\alpha_j:=1-\displaystyle\frac{1}{p_j}$. If $1<p_j<\infty$, the result just proved with second index $p_j$, together with the final identity in Lemma~\ref{lem:identidades-reordenamiento-decreciente}, identifies
\[
L^{p_j}(A)=L^{p_j,p_j}(A)
=(L^1(A),L^\infty(A))_{\alpha_j,p_j}
\]
with equivalent norms. If $p_0=1$ or $p_1=\infty$, retain the corresponding endpoint. Apply Theorem~\ref{teo:reiteracion-interpolacion-real}; its final parameter is
\[
(1-\theta)\alpha_0+\theta\alpha_1
=1-\left(\frac{1-\theta}{p_0}+\frac{\theta}{p_1}\right)=1-\frac{1}{p}.
\]
We obtain the asserted identity. Equivalence of norms at the endpoints is transported by Theorem~\ref{teo:interpolacion-operadores-metodo-K} applied to the identity in both directions; this gives the indicated dependence of the constants.

The equality $L^{p,p}(A)=L^p(A)$ is isometric for the quasinorm in Definition~\ref{def:espacios-lorentz}, since
\[
\|f\|_{L^{p,p}(A)}^p=\int_0^\infty f^*(t)^p\,dt=\|f\|_{L^p(A)}^p.
\]
Finally, monotonicity of $f^*$ allows comparison, on each interval $[2^k,2^{k+1}]$, of the norm defining $L^{p,q}(A)$ with that of the sequence $(2^{\displaystyle\frac{k}{p}}f^*(2^k))_k$ in $\ell^q(\mathbb Z)$. The lower bound uses the right endpoint and the upper bound the left; shifting the index only introduces the factor $2^{\displaystyle\frac{1}{p}}$. The embedding $\ell^{q_0}(\mathbb Z)\hookrightarrow\ell^{q_1}(\mathbb Z)$ when $q_0\leq q_1$ proves the final assertion.
\end{proof}

\chapter{Intermediate regularity}
\label{cap:regularidad-intermedia}

Regularity does not increase only in integer steps. Between having $m$ and $m+1$ integrable derivatives lies a continuous range of behaviors arising in trace problems, interpolation, nonlocal equations, and spectral theory. To measure these intermediate orders, we shall use differences of values, Fourier multipliers, and frequency decompositions.

\begin{semblanzaHistorica}{Several scales for the same question}
Bessel potentials connected regularity with the Fourier transform; Slobodeckij and Gagliardo developed difference-based seminorms; Besov and Triebel organized scales capable of separating the order of smoothness from the way frequencies are summed. Each construction reflects different properties of functions. Their agreements and differences allow us to choose the appropriate scale for localization in charts, taking traces, and interpolating operators.
\end{semblanzaHistorica}

The spaces $W^{m,p}(\mathbb R^n)$ measure regularity through integer-order derivatives. Nonlocal problems show that intermediate orders also need to be quantified. For $0<s<1$, the formal equation
\[
(-\Delta)^s u=f
\]
is associated with the energy
\[
\mathcal E_s(u)
=
\frac{c_{n,s}}{4}
\int_{\mathbb R^n}\int_{\mathbb R^n}
\frac{|u(x)-u(y)|^2}{|x-y|^{n+2s}}\,dx\,dy
-\int_{\mathbb R^n}fu\,dx.
\]
The first term compares values of $u$ at pairs of points instead of using a local derivative, leading to the Gagliardo--Slobodeckij seminorm. For $1<p<\infty$, the analogous expression
\[
\int_{\mathbb R^n}\int_{\mathbb R^n}
\frac{|u(x)-u(y)|^p}{|x-y|^{n+sp}}\,dx\,dy
\]
motivates the spaces $W^{s,p}(\mathbb R^n)$ and appears in energies associated with the fractional $p$-Laplacian.

There is another route to fractional regularity. The operator $(I-\Delta)^{\displaystyle\frac{s}{2}}$ is defined through Fourier analysis and leads to Bessel potentials and the spaces $H^{s,p}(\mathbb R^n)$. Both constructions recover $W^{m,p}(\mathbb R^n)$ when $s=m\in\mathbb N_0$, but for $s\notin\mathbb N_0$ they agree in general only in the Hilbert case $p=2$. Traces and elliptic regularity also require negative orders and scales well suited to localization.

This is why Sobolev--Slobodeckij, Bessel, Besov, and Triebel--Lizorkin spaces appear together. Each scale organizes integrability, smoothness order, and frequency distribution differently. This chapter constructs these scales, relates their parameters, and proves the interpolation, localization, and embedding properties to be used on manifolds. We follow \cite[Chapters~6--8 and 11--12]{Leoni2023}, \cite[Chapter~1]{Triebel2}, and \cite[\S\S2.2--2.3]{Schneider2021}.

\section{Bessel potentials and the spaces \texorpdfstring{$H^{s,p}(\mathbb R^n)$}{Hsp}}
\label{sec:potenciales-bessel-euclidianos}

Write
\[
\langle\xi\rangle:=(1+\|\xi\|^2)^{\frac{1}{2}},
\qquad \xi\in\mathbb R^n.
\]

\begin{lemma}[Peetre's inequality]
\label{lem:desigualdad-peetre-peso-japones}
\index{Peetre inequality@Peetre's inequality}
For every $\mu\in\mathbb R$ and any $\xi,\eta\in\mathbb R^n$,
\[
\langle\xi+\eta\rangle^\mu
\leq
2^{\frac{|\mu|}{2}}
\langle\xi\rangle^\mu
\langle\eta\rangle^{|\mu|}.
\]
\end{lemma}

\begin{proof}
The triangle inequality and $2ab\leq a^2+b^2$ give
\[
\begin{aligned}
\langle\xi+\eta\rangle^2
&=1+\|\xi+\eta\|^2\\
&\leq 1+2\|\xi\|^2+2\|\eta\|^2\\
&\leq 2(1+\|\xi\|^2)(1+\|\eta\|^2)
=2\langle\xi\rangle^2\langle\eta\rangle^2.
\end{aligned}
\]
If $\mu\geq0$, raise this inequality to the power $\displaystyle\frac{\mu}{2}$. If $\mu<0$, apply the preceding inequality to $\xi=(\xi+\eta)-\eta$ to obtain
\[
\langle\xi\rangle
\leq
\sqrt2\,\langle\xi+\eta\rangle\langle\eta\rangle.
\]
Dividing by $\langle\xi\rangle\langle\xi+\eta\rangle$ and raising to the power $|\mu|$ gives
\[
\langle\xi+\eta\rangle^\mu
\leq
2^{\frac{|\mu|}{2}}
\langle\xi\rangle^\mu
\langle\eta\rangle^{|\mu|},
\]
as required.
\end{proof}

For each $s\in\mathbb R$ and each multi-index $\alpha$, there is $C_{s,\alpha}>0$ such that
\[
|D^\alpha\langle\xi\rangle^s|
\leq C_{s,\alpha}\langle\xi\rangle^{s-|\alpha|},
\qquad \xi\in\mathbb R^n.
\]
In particular, all these derivatives have at most polynomial growth. By the remark following Definition~\ref{def:operaciones-distribuciones-temperadas}, multiplication $T\mapsto\langle\,\cdot\,\rangle^sT$ is well defined and continuous on $\mathcal S'(\mathbb R^n)$.

\begin{definition}[Bessel operator and space]
\label{def:regularidad-intermedia-norma-espacio-sobolev-fraccionario}
\glsadd{espacio-bessel}
\index{Bessel potential}
\index{Bessel potential space}
Let $s\in\mathbb R$ and $1<p<\infty$. Define
\[
J^su
:=
\mathcal F^{-1}\bigl(\langle\,\cdot\,\rangle^s\widehat u\bigr),
\qquad u\in\mathcal S'(\mathbb R^n),
\]
and
\[
H^{s,p}(\mathbb R^n)
:=
\{u\in\mathcal S'(\mathbb R^n)\mid J^su\in L^p(\mathbb R^n)\},
\qquad
\|u\|_{H^{s,p}(\mathbb R^n)}
:=
\|J^su\|_{L^p(\mathbb R^n)}.
\]
The operator $J^{-s}$ is called the \textit{Bessel potential} of order $s$.
\end{definition}

\begin{proposition}[Structure of the Bessel scale]
\label{prop:estructura-escala-bessel-euclidiana}
Let $s,t\in\mathbb R$ and $1<p<\infty$.
\begin{enumerate}[label=(\alph*)]
\item $H^{s,p}(\mathbb R^n)$ is a Banach space, and $J^s\colon H^{s,p}(\mathbb R^n)\longrightarrow L^p(\mathbb R^n)$ is an isometric isomorphism with inverse $J^{-s}$.
\item $J^t\colon H^{s,p}(\mathbb R^n)\longrightarrow H^{s-t,p}(\mathbb R^n)$ is an isometric isomorphism.
\item $\mathcal S(\mathbb R^n)$ is dense in $H^{s,p}(\mathbb R^n)$.
\item With respect to the distributional pairing extended from $\mathcal S(\mathbb R^n)$, the dual of $H^{s,p}(\mathbb R^n)$ is identified isometrically with $H^{-s,p'}(\mathbb R^n)$: if $\Lambda_v(u):=\langle v,u\rangle$, then
\[
 \|\Lambda_v\|_{(H^{s,p}(\mathbb R^n))'}=\|v\|_{H^{-s,p'}(\mathbb R^n)}.
\]
\end{enumerate}
\end{proposition}

\begin{proof}
The symbols $\langle\xi\rangle^s$ and $\langle\xi\rangle^{-s}$
act continuously by multiplication on $\mathcal S$ and $\mathcal S'$,
and their product is one. Thus $J^s$ and $J^{-s}$ are inverses.
By definition, $J^s\colon H^{s,p}\to L^p$ is a surjective isometry
with inverse $J^{-s}$. Completeness of $L^p$ proves (a).
For (b), simply observe that
\[
J^{s-t}J^t u=J^s u,
\]
so $J^t$ maps $H^{s,p}$ isometrically onto $H^{s-t,p}$.
For (c), if $f_k\in\mathcal S$ converges to $J^s u$ in $L^p$,
then $J^{-s}f_k\in\mathcal S$ and converges to $u$ in $H^{s,p}$.

We make the pairing in (d) explicit. The symbol of $J^s$ is even;
the Fourier duality formulas imply
\[
\langle J^s T,\phi\rangle=\langle T,J^s\phi\rangle,
\qquad T\in\mathcal S',\quad\phi\in\mathcal S.
\]
Indeed, the transpose of $\mathcal F^{-1}a\mathcal F$
has symbol $a(-\xi)$, which here equals $a(\xi)$.
If $v\in H^{-s,p'}$ and $g=J^{-s}v\in L^{p'}$, then
$v=J^s g$ and
\[
\langle v,\phi\rangle=\int_{\mathbb R^n}g(x)J^s\phi(x)\,dx.
\]
Hölder's inequality extends this expression to a functional on
$H^{s,p}$ with norm at most $\|g\|_{p'}$. Conversely, a functional
$\Lambda\in(H^{s,p})'$ induces
$f\mapsto\Lambda(J^{-s}f)$ on $L^p$, with the same norm.
Duality of $L^p$ gives a unique $g\in L^{p'}$ of that norm.
The distribution $v=J^s g$ belongs to $H^{-s,p'}$ and represents
$\Lambda$ by the preceding formula, first on $\mathcal S$ and
then on all of $H^{s,p}$ by (c). This construction proves both
isometry and uniqueness in (d).
\end{proof}

For $s>0$, the potential $J^{-s}$ admits a convolution representation. This description links the Fourier construction to the later construction through the heat semigroup.

\begin{proposition}[Euclidean Bessel kernel]
\label{prop:nucleo-bessel-euclidiano}
Let $s>0$. There is a positive radial function $G_s\in L^1(\mathbb R^n)$, with $\displaystyle\int_{\mathbb R^n}G_s=1$, such that
\[
J^{-s}f=G_s*f
\]
for $f\in\mathcal S(\mathbb R^n)$ and, by continuity, for $f\in L^p(\mathbb R^n)$, $1\leq p\leq\infty$. Moreover,
\[
G_s(x)
=
\frac{1}{\Gamma(\frac{s}{2})}
\int_0^\infty
e^{-t}(4\pi t)^{-\frac{n}{2}}e^{-\frac{\|x\|^2}{4t}}
t^{\frac{s}{2}-1}\,dt,
\qquad x\in\mathbb R^n\setminus\{0\}.
\]

The value at the origin may be chosen arbitrarily without changing the element of \(L^1(\mathbb R^n)\); the integral there may diverge when \(s\leq n\).
\end{proposition}

\begin{proof}
Set $a=\displaystyle\frac{s}{2}>0$ and
\[
H_t(x):=(4\pi t)^{-\frac{n}{2}}e^{-\frac{\|x\|^2}{4t}},\qquad
G_s(x):=\frac{1}{\Gamma(a)}\int_0^\infty t^{a-1}e^{-t}H_t(x)\,dt.
\]
Each $H_t$ is nonnegative and has integral one.
Tonelli's theorem gives
\[
\int_{\mathbb R^n}G_s(x)\,dx
=\frac{1}{\Gamma(a)}\int_0^\infty t^{a-1}e^{-t}\,dt=1.
\]
Thus the formula defines an $L^1$ function, finite almost everywhere,
radial and nonnegative. For $x\neq0$, it is also finite pointwise:
$e^{-\displaystyle\frac{\|x\|^2}{4t}}$ controls $t\to0^+$ and $e^{-t}$ controls
$t\to\infty$.

With the unitary normalization,
$\widehat H_t(\xi)=(2\pi)^{-\displaystyle\frac{n}{2}}e^{-t\|\xi\|^2}$.
The joint absolute integral is finite by the mass calculation,
so Fubini's theorem and $r=t(1+\|\xi\|^2)$ give
\[
\widehat G_s(\xi)
=\frac{(2\pi)^{-\frac{n}{2}}}{\Gamma(a)}
\int_0^\infty t^{a-1}e^{-t(1+\|\xi\|^2)}\,dt
=(2\pi)^{-\frac{n}{2}}\langle\xi\rangle^{-s}.
\]
Moreover, for every $N\in\mathbb N_0$, the change $x=\sqrt t\,y$
in the heat kernel shows that
\[
\int_{\mathbb R^n}(1+\|x\|)^N H_t(x)\,dx
\leq C_N(1+t^{\frac{N}{2}}).
\]
Integrating this bound against $t^{a-1}\displaystyle\frac{e^{-t}}{\Gamma(a)}$ proves
that $G_s$ has finite absolute moments of every polynomial order.
Hence, for $\psi\in\mathcal S$, differentiation of $G_s*\psi$
under the integral and the weight inequality give
\[
q_N(G_s*\psi)
\leq C_N q_N(\psi)
\int_{\mathbb R^n}(1+\|y\|)^N G_s(y)\,dy.
\]
In particular, $G_s*\psi\in\mathcal S$. The Fourier convolution
formula, justified by Fubini for the $L^1$ functions involved,
gives $G_s*\psi=J^{-s}\psi$.

For $1\leq p<\infty$, Young's inequality gives
$\|G_s*f\|_p\leq\|f\|_p$. Approximate $f\in L^p$ by Schwartz
functions. Their convolutions converge in $L^p$ by Young's inequality
and in $\mathcal S'$ to $J^{-s}f$ by continuity of $J^{-s}$.
The continuous inclusion $L^p\hookrightarrow\mathcal S'$ identifies
the limits.
If $p=\infty$, Fubini applies for $\psi\in\mathcal S$ because
the absolute integral is bounded by
$\|f\|_\infty\|G_s\|_1\|\psi\|_1$. Since $G_s$ is even,
\[
\langle G_s*f,\psi\rangle
=\int_{\mathbb R^n}f(x)(G_s*\psi)(x)\,dx
=\langle f,J^{-s}\psi\rangle
=\langle J^{-s}f,\psi\rangle.
\]
The last identity is the self-transpose property proved in
Proposition~\ref{prop:estructura-escala-bessel-euclidiana}.
Young's inequality also gives the norm bound for $p=\infty$.
\end{proof}

\begin{definition}[Bessel spaces on an open subset]
\label{def:regularidad-intermedia-norma-r-n-abierto-s}
\index{Bessel potential space!on an open subset}
Let $\Omega\subseteq\mathbb R^n$ be open. Define $H^{s,p}(\Omega)$ as the space of restrictions to $\Omega$ of elements of $H^{s,p}(\mathbb R^n)$, with the quotient norm
\[
\|u\|_{H^{s,p}(\Omega)}
:=
\inf\left\{
\|U\|_{H^{s,p}(\mathbb R^n)}
\;\middle|\;
U\restriction_\Omega=u
\right\}.
\]
\end{definition}

This definition is extrinsic: it does not assert that a seminorm computed only within $\Omega$ is equivalent to the quotient norm. Such an equivalence requires an extension property of the open subset and will be proved for Lipschitz domains. This distinction will be essential when comparing $H^{s,p}(\Omega)$ with intrinsic Slobodeckij spaces.

\subsection{The Euclidean Hilbert scale and localization}

For elliptic regularity, we shall need only the case $p=2$. Write $H^s(\mathbb R^n):=H^{s,2}(\mathbb R^n)$ and $\langle\xi\rangle:=(1+\|\xi\|^2)^{\displaystyle\frac{1}{2}}$. If $r\in\mathbb N$, define the space $H^s(\mathbb R^n,\mathbb K^r)$ componentwise and equip it with the norm
\[
\|u\|_{H^s(\mathbb R^n,\mathbb K^r)}^2
:=
\displaystyle\sum_{a=1}^r
\|u^a\|_{H^s(\mathbb R^n)}^2.
\]
Although the preceding definition is available for every $s\in\mathbb R$, the proof of local elliptic regularity will use only integer indices, including negative ones. Noninteger orders will be revisited when studying the scales associated with the heat semigroup.

\begin{proposition}[Elementary properties of the $H^s(\mathbb R^n)$ scale]
\label{prop:regularidad-intermedia-propiedades-escala-hilbertiana}
Let $s,t\in\mathbb R$ and $t\geq s$, and let $\alpha$ be a multi-index. Then the following assertions hold.
\begin{enumerate}[label=(\alph*)]
\item $H^s(\mathbb R^n,\mathbb K^r)$ is a Hilbert space, and $\mathcal S(\mathbb R^n,\mathbb K^r)$ is dense in it.
\item The inclusion $H^t(\mathbb R^n,\mathbb K^r)\hookrightarrow H^s(\mathbb R^n,\mathbb K^r)$ is continuous.
\item The distributional derivative induces a continuous linear operator
\[
D^\alpha:H^s(\mathbb R^n,\mathbb K^r)
\longrightarrow
H^{s-|\alpha|}(\mathbb R^n,\mathbb K^r).
\]
\item Translation $\tau_hu(x):=u(x-h)$ is an isometry of $H^s(\mathbb R^n,\mathbb K^r)$, and
\[
\|\tau_hu-u\|_{H^s(\mathbb R^n,\mathbb K^r)}\longrightarrow0
\qquad\text{when }h\to0.
\]
\item In the real case, the dual of $H^s(\mathbb R^n,\mathbb R^r)$ is identified isometrically with $H^{-s}(\mathbb R^n,\mathbb R^r)$. In the complex case, the correct complex-linear Hermitian identification is
\begin{equation}
\label{eq:dualidad-hermitiana-bessel-euclidiana}
 H^{-s}(\mathbb R^n,\overline{\mathbb C^r})
 \xrightarrow{\ \cong\ }
 \bigl(H^s(\mathbb R^n,\mathbb C^r)\bigr)'.
\end{equation}
Moreover, the standard basis of $\mathbb C^r$ defines the complex-bilinear distributional pairing
\begin{equation}
\label{eq:dualidad-bilineal-bessel-euclidiana}
 H^{-s}(\mathbb R^n,\mathbb C^r)
 \times H^s(\mathbb R^n,\mathbb C^r)\longrightarrow\mathbb C,
 \qquad (v,u)\longmapsto\langle v,u\rangle_{\mathcal D',\mathcal D},
\end{equation}
which also induces an isometric complex-linear isomorphism from the first factor onto the dual of the second. This second identification uses the standard trivialization; the first is the intrinsic Hermitian form.
\end{enumerate}
\end{proposition}

\begin{proof}
By Plancherel's theorem~\ref{teo:plancherel-fourier}, the map
\[
\Lambda^s:H^s(\mathbb R^n,\mathbb K^r)
\longrightarrow L^2(\mathbb R^n,\mathbb K^r),
\qquad
\widehat{\Lambda^su}(\xi):=\langle\xi\rangle^s\widehat u(\xi),
\]
is a surjective isometry with inverse $\Lambda^{-s}$. In the real case, complexification is used when taking Fourier transforms. The multipliers $\langle\xi\rangle^{\pm s}$, being real and even, commute with conjugation in physical space and preserve real distributions. Thus the isometry and its inverse are well defined also for $\mathbb K=\mathbb R$. Completeness of $L^2(\mathbb R^n,\mathbb K^r)$ proves that $H^s(\mathbb R^n,\mathbb K^r)$ is a Hilbert space with the inner product transported by $\Lambda^s$.

For density, let $u\in H^s(\mathbb R^n,\mathbb K^r)$. Density in $L^2(\mathbb R^n,\mathbb K^r)$ allows us to choose $(f_j)\subseteq C_c^\infty(\mathbb R^n,\mathbb K^r)$ such that $f_j\to\Lambda^su$ in $L^2(\mathbb R^n,\mathbb K^r)$. Set $u_j:=\Lambda^{-s}f_j$. Derivatives of the symbol $\langle\xi\rangle^{-s}$ have at most polynomial growth, so multiplication by it preserves the Schwartz class. Since the Fourier transform and its inverse also preserve that class, $u_j\in\mathcal S(\mathbb R^n,\mathbb C^r)$; if $\mathbb K=\mathbb R$, compatibility with conjugation also gives $u_j=\overline{u_j}$. Thus, in both cases, $u_j\in\mathcal S(\mathbb R^n,\mathbb K^r)$ and
\[
\|u_j-u\|_{H^s(\mathbb R^n,\mathbb K^r)}
=
\|f_j-\Lambda^su\|_{L^2(\mathbb R^n,\mathbb K^r)}
\longrightarrow0.
\]

If $t\geq s$, then $\langle\xi\rangle^s\leq\langle\xi\rangle^t$, whence
\[
\|u\|_{H^s(\mathbb R^n,\mathbb K^r)}
=
\|\langle\xi\rangle^s\widehat u\|_{L^2(\mathbb R^n,\mathbb C^r)}
\leq
\|\langle\xi\rangle^t\widehat u\|_{L^2(\mathbb R^n,\mathbb C^r)}
=
\|u\|_{H^t(\mathbb R^n,\mathbb K^r)}.
\]
Likewise, since $\widehat{D^\alpha u}(\xi)=(i\xi)^\alpha\widehat u(\xi)$ and $|\xi^\alpha|\leq\langle\xi\rangle^{|\alpha|}$,
\[
\|D^\alpha u\|_{H^{s-|\alpha|}(\mathbb R^n,\mathbb K^r)}
\leq
\|u\|_{H^s(\mathbb R^n,\mathbb K^r)}.
\]

The Fourier transform of $\tau_hu$ is $e^{-ih\cdot\xi}\widehat u(\xi)$, so $\|\tau_hu\|_{H^s(\mathbb R^n,\mathbb K^r)}
=\|u\|_{H^s(\mathbb R^n,\mathbb K^r)}$. Moreover,
\[
\|\tau_hu-u\|_{H^s(\mathbb R^n,\mathbb K^r)}^2
=
\int_{\mathbb R^n}
\langle\xi\rangle^{2s}|e^{-ih\cdot\xi}-1|^2|\widehat u(\xi)|^2\,d\xi.
\]
The integrand converges pointwise to zero and is dominated by $4\langle\xi\rangle^{2s}\|\widehat u(\xi)\|_{\mathbb C^r}^2
\in L^1(\mathbb R^n)$. Lebesgue's dominated convergence theorem~\ref{convergencia dominada} gives continuity of translations.

Finally, first let $\mathbb K=\mathbb C$. By the complex-linear form of Riesz's theorem, Corollary~\ref{cor:riesz-lineal-desde-conjugado}, every functional $F$ on $H^s(\mathbb R^n,\mathbb C^r)$ has a unique representative $\overline g\in\overline{L^2(\mathbb R^n,\mathbb C^r)}$ such that
\[
 F(u)=\langle\Lambda^su,g\rangle_{L^2(\mathbb R^n,\mathbb C^r)}.
\]
If $w:=\Lambda^sg\in H^{-s}(\mathbb R^n,\mathbb C^r)$, this equality reads
\[
 F(u)=\langle\Lambda^su,\Lambda^{-s}w\rangle_{L^2(\mathbb R^n,\mathbb C^r)}.
\]
Pointwise conjugation identifies $\overline{H^{-s}(\mathbb R^n,\mathbb C^r)}$ complex-linearly with $H^{-s}(\mathbb R^n,\overline{\mathbb C^r})$; this gives \eqref{eq:dualidad-hermitiana-bessel-euclidiana}, including equality of norms.

For the bilinear pairing, if $u,v$ are Schwartz functions, Plancherel's theorem and the Fourier convention give
\[
 \langle v,u\rangle_{\mathcal D',\mathcal D}
 =\int_{\mathbb R^n}\sum_{a=1}^r
   \widehat v^{,a}(-\xi)\widehat u^{,a}(\xi)\,d\xi.
\]
Inserting the weights $\langle\xi\rangle^{-s}$ and $\langle\xi\rangle^s$, Cauchy--Schwarz and the substitution $\xi\mapsto-\xi$ prove that its norm is at most $\|v\|_{H^{-s}(\mathbb R^n,\mathbb C^r)}\|u\|_{H^s(\mathbb R^n,\mathbb C^r)}$. Equality of norms and surjectivity follow from the preceding Hermitian identification through conjugation and reflection $\xi\mapsto-\xi$. Density of $\mathcal S$ extends the pairing and proves \eqref{eq:dualidad-bilineal-bessel-euclidiana}. In the real case, omit conjugation, and the same argument gives the stated assertion.
\end{proof}

\begin{definition}[Local Sobolev spaces of real order]
\label{def:regularidad-intermedia-sobolev-local-orden-real}
Let $\Omega\subseteq\mathbb R^n$ be open and let $s\in\mathbb R$. Define
\[
H^s_{\mathrm{loc}}(\Omega,\mathbb K^r)
:=
\left\{
 u\in\mathcal D'(\Omega,\mathbb K^r)
 \;\middle|\;
 \chi u\in H^s(\mathbb R^n,\mathbb K^r)
 \text{ for every }\chi\in C_c^\infty(\Omega)
\right\},
\]
where $\chi u$ is extended by zero outside $\Omega$.
\end{definition}

\subsection{Multiplication, regularization, and commutators}

Ordinary continuity of multiplication by a smooth function does not suffice to freeze coefficients at negative orders. We shall need to separate the principal term, controlled by the uniform norm of the coefficient, from a remainder losing one derivative.

\begin{lemma}[Refined multiplication estimate]
\label{lem:regularidad-intermedia-multiplicacion-afinada}
Let $s\in\mathbb R$ and $a\in C_c^\infty(\mathbb R^n)$. There is a constant $C_{a,s}>0$ such that
\[
\|av\|_{H^s(\mathbb R^n)}
\leq
\|a\|_{L^\infty(\mathbb R^n)}\|v\|_{H^s(\mathbb R^n)}
+C_{a,s}\|v\|_{H^{s-1}(\mathbb R^n)}
\]
for every $v\in H^s(\mathbb R^n)$. For a matrix-valued function $A\in C_c^\infty(\mathbb R^n,M_{q\times r}(\mathbb K))$, there is a constant $C_{A,s}>0$ such that
\[
\|Av\|_{H^s(\mathbb R^n,\mathbb K^q)}
\leq
C_{q,r}\|A\|_{L^\infty(\mathbb R^n,M_{q\times r}(\mathbb K))}\|v\|_{H^s(\mathbb R^n,\mathbb K^r)}
+C_{A,s}\|v\|_{H^{s-1}(\mathbb R^n,\mathbb K^r)}.
\]
\end{lemma}

\begin{proof}
By density, it suffices to consider $v\in\mathcal S(\mathbb R^n)$. Write $\Lambda^s$ for the Fourier multiplier with symbol $\langle\xi\rangle^s$. Then
\[
\Lambda^s(av)=a\Lambda^sv+R_s(a,v),
\]
where, up to the fixed constant determined by the normalization of the Fourier transform,
\[
\widehat{R_s(a,v)}(\xi)
=
\int_{\mathbb R^n}
\bigl(\langle\xi\rangle^s-\langle\xi-\eta\rangle^s\bigr)
\widehat a(\eta)\widehat v(\xi-\eta)\,d\eta.
\]
The fundamental theorem of calculus for Banach-valued functions~\ref{teo:b5-fundamental-calculo-riemann-banach}, applied to $t\mapsto\langle\xi-t\eta\rangle^s$, together with Peetre's inequality from Lemma~\ref{lem:desigualdad-peetre-peso-japones}, provides a constant $C_s>0$ for which
\[
\left|\langle\xi\rangle^s-\langle\xi-\eta\rangle^s\right|
\leq
C_s\|\eta\|\langle\eta\rangle^{|s|+1}
\langle\xi-\eta\rangle^{s-1}.
\]
Since $a\in C_c^\infty$, its Fourier transform is rapidly decreasing, and the function
\[
K_s(\eta):=C_s\|\eta\|\langle\eta\rangle^{|s|+1}|\widehat a(\eta)|
\]
belongs to $L^1(\mathbb R^n)$. Young's convolution inequality from Theorem~\ref{teo:young-convoluciones} gives
\[
\|R_s(a,v)\|_{L^2(\mathbb R^n)}
\leq
\|K_s\|_{L^1(\mathbb R^n)}\|v\|_{H^{s-1}(\mathbb R^n)}.
\]
On the other hand,
\[
\|a\Lambda^sv\|_{L^2(\mathbb R^n)}
\leq
\|a\|_{L^\infty(\mathbb R^n)}\|v\|_{H^s(\mathbb R^n)}.
\]
Adding the two estimates proves the scalar inequality. The matrix assertion follows by applying the result to each entry and summing finitely many components.
\end{proof}

\begin{lemma}[Localization criterion]
\label{lem:regularidad-intermedia-criterio-localizacion-sobolev}
Let $u\in\mathcal D'(\Omega,\mathbb K^r)$ and $s\in\mathbb R$. The following conditions are equivalent.
\begin{enumerate}[label=(\alph*)]
\item $u\in H^s_{\mathrm{loc}}(\Omega,\mathbb K^r)$.
\item For every compact set $K\subseteq\Omega$, there is $\chi\in C_c^\infty(\Omega)$, with $\chi=1$ on a neighborhood of $K$, such that $\chi u\in H^s(\mathbb R^n,\mathbb K^r)$.
\item Every point $x\in\Omega$ has an open neighborhood $V\Subset\Omega$ and a function $\chi\in C_c^\infty(\Omega)$, with $\chi=1$ on a neighborhood of $\overline V$, for which $\chi u\in H^s(\mathbb R^n,\mathbb K^r)$.
\end{enumerate}
\end{lemma}

\begin{proof}
The implication (a)$\Rightarrow$(b) is immediate after choosing, by Proposition~\ref{funciones flan}, a bump function equal to one on a neighborhood of $K$; (b)$\Rightarrow$(c) follows by taking a compact set containing the point in its interior.

Assume (c) and let $\zeta\in C_c^\infty(\Omega)$. The compact set $\operatorname{supp}(\zeta)$ is covered by finitely many open sets $V_1,\dots,V_N$ as in (c). Choose a partition of unity $\eta_1,\dots,\eta_N$ on a neighborhood of $\operatorname{supp}(\zeta)$ subordinate to this cover, and functions $\chi_j\in C_c^\infty(\Omega)$ such that $\chi_j=1$ on a neighborhood of $\overline{V_j}$ and $\chi_ju\in H^s(\mathbb R^n,\mathbb K^r)$. Then
\[
\zeta u
=
\displaystyle\sum_{j=1}^N(\zeta\eta_j)(\chi_ju).
\]
The refined multiplication estimate in Lemma~\ref{lem:regularidad-intermedia-multiplicacion-afinada} shows that each summand belongs to $H^s(\mathbb R^n,\mathbb K^r)$. The sum is finite, so $\zeta u\in H^s(\mathbb R^n,\mathbb K^r)$. Therefore, $u\in H^s_{\mathrm{loc}}(\Omega,\mathbb K^r)$.
\end{proof}

\begin{proposition}[Local Sobolev order of a distribution]
\label{prop:regularidad-intermedia-distribucion-orden-sobolev-negativo}
Let $\Omega\subseteq\mathbb R^n$ be open and let $T\in\mathcal D'(\Omega,\mathbb K^r)$. For every $\chi\in C_c^\infty(\Omega)$, there is $N\in\mathbb N$ such that
\[
\chi T\in H^{-N}(\mathbb R^n,\mathbb K^r).
\]
In particular, every distribution belongs locally to some Sobolev space of negative integer order, although the integer may depend on the compact set under consideration.
\end{proposition}

\begin{proof}
The argument is componentwise, so it suffices to treat the scalar case. Let $K\Subset\Omega$ be a compact set containing the support of $\chi$ in its interior. Continuity of $T$ on $\mathcal D_K(\Omega)$ provides an integer $q\in\mathbb N_0$ and a constant $C>0$ such that
\[
|\langle T,\psi\rangle|
\leq
C\max_{\substack{\alpha\in\mathbb N_0^n\\|\alpha|\leq q}}\sup_{x\in K}|D^\alpha\psi(x)|
\qquad
\text{for every }\psi\in\mathcal D_K(\Omega).
\]
Applying this inequality to $\psi=\chi\varphi$, with $\varphi\in \mathcal S(\mathbb R^n)$, and using the Leibniz rule gives $C_1>0$ such that
\[
|\langle\chi T,\varphi\rangle|
\leq
C_1\max_{\substack{\beta\in\mathbb N_0^n\\|\beta|\leq q}}\|D^\beta\varphi\|_{L^\infty(\mathbb R^n)}.
\]
Choose an integer $N>q+\displaystyle\frac{n}{2}$. By the Fourier inversion theorem~\ref{teo:inversion-fourier-schwartz} and the Cauchy--Schwarz inequality from Theorem~\ref{teo:cauchy-schwarz-hilbert}, for $|\beta|\leq q$ we have
\[
|D^\beta\varphi(x)|
\leq
C_2\int_{\mathbb R^n}\|\xi\|^{|\beta|}|\widehat\varphi(\xi)|\,d\xi
\leq
C_3
\left(
\int_{\mathbb R^n}\langle\xi\rangle^{-2(N-|\beta|)}\,d\xi
\right)^{\frac{1}{2}}
\|\varphi\|_{H^N(\mathbb R^n)}.
\]
The integral is finite because $N-|\beta|>\displaystyle\frac{n}{2}$. Thus,
\[
|\langle\chi T,\varphi\rangle|
\leq
C_4\|\varphi\|_{H^N(\mathbb R^n)}.
\]
Density of $\mathcal S(\mathbb R^n)$ in $H^N(\mathbb R^n)$, proved in Proposition~\ref{prop:regularidad-intermedia-propiedades-escala-hilbertiana}, extends $\chi T$ to a continuous functional on $H^N(\mathbb R^n)$. The dual identification in Proposition~\ref{prop:regularidad-intermedia-propiedades-escala-hilbertiana} yields $\chi T\in H^{-N}(\mathbb R^n)$.
\end{proof}

Let $\rho$ be the radial mollifier from Theorem~\ref{teo:molificadores-euclidianos}, and set $J_\varepsilon u:=\rho_\varepsilon*u$.

\begin{proposition}[Mollifiers on the $H^s(\mathbb R^n)$ scale]
\label{prop:regularidad-intermedia-molificadores-escala-sobolev}
Let $s\in\mathbb R$ and $0<\varepsilon\leq1$. Then $J_\varepsilon$ induces a continuous linear operator on $H^s(\mathbb R^n,\mathbb K^r)$, and
\[
\|J_\varepsilon u\|_{H^s(\mathbb R^n,\mathbb K^r)}
\leq
\|u\|_{H^s(\mathbb R^n,\mathbb K^r)}.
\]
Moreover,
\[
J_\varepsilon u\longrightarrow u
\quad\text{in }H^s(\mathbb R^n,\mathbb K^r)
\quad\text{when }\varepsilon\to0^+.
\]
If $u$ has compact support, then $J_\varepsilon u\in C^\infty(\mathbb R^n,\mathbb K^r)$ and
\[
\operatorname{supp}(J_\varepsilon u)
\subseteq
\left\{x\in\mathbb R^n\mid
\operatorname{dist}_{\mathrm{euc}}(x,\operatorname{supp}(u))\leq\varepsilon
\right\}.
\]
Finally, since $\rho$ is radial, $J_\varepsilon$ is self-adjoint with respect to the $L^2(\mathbb R^n,\mathbb K^r)$ pairing and the $H^{-s}(\mathbb R^n,\mathbb K^r)$--$H^s(\mathbb R^n,\mathbb K^r)$ duality.
\end{proposition}

\begin{proof}
The change of variable $y=\varepsilon z$ allows us to write
\[
J_\varepsilon u
=
\int_{\mathbb R^n}\rho(z)\tau_{\varepsilon z}u\,dz,
\]
where the integral is taken in the Hilbert space $H^s(\mathbb R^n,\mathbb K^r)$. Since translations are isometries,
\[
\|J_\varepsilon u\|_{H^s(\mathbb R^n,\mathbb K^r)}
\leq
\int_{\mathbb R^n}\rho(z)
\|\tau_{\varepsilon z}u\|_{H^s(\mathbb R^n,\mathbb K^r)}\,dz
=
\|u\|_{H^s(\mathbb R^n,\mathbb K^r)}.
\]
Likewise,
\[
\|J_\varepsilon u-u\|_{H^s(\mathbb R^n,\mathbb K^r)}
\leq
\int_{\mathbb R^n}\rho(z)
\|\tau_{\varepsilon z}u-u\|_{H^s(\mathbb R^n,\mathbb K^r)}\,dz.
\]
For each fixed $z$, the integrand converges to zero by Proposition~\ref{prop:regularidad-intermedia-propiedades-escala-hilbertiana} and is dominated by $2\rho(z)\|u\|_{H^s(\mathbb R^n,\mathbb K^r)}$. The Bochner dominated convergence theorem~\ref{teo:b5-convergencia-dominada-bochner} proves convergence in $H^s(\mathbb R^n,\mathbb K^r)$.

If $u$ has compact support, convolution of a compactly supported distribution with $\rho_\varepsilon\in C_c^\infty(\mathbb R^n)$ is smooth. The support inclusion follows directly from the fact that $\rho_\varepsilon$ vanishes outside $\overline B_{\mathrm{euc}}(0,\varepsilon)$. Finally, radiality implies $\rho_\varepsilon(x)=\rho_\varepsilon(-x)$; Fubini's theorem~\ref{Fubini}, applied first to smooth functions, and density give $\langle J_\varepsilon u,v\rangle=\langle u,J_\varepsilon v\rangle$ and then the same identity by duality.
\end{proof}

\begin{lemma}[Friedrichs commutator]
\label{lem:regularidad-intermedia-conmutador-friedrichs}
Let $j\in\mathbb Z$ and $a\in C_c^\infty(\mathbb R^n)$. There is a constant $C_{a,j}>0$, independent of $0<\varepsilon\leq1$, such that
\[
\|[a,J_\varepsilon]v\|_{H^{j+1}(\mathbb R^n)}
\leq
C_{a,j}\|v\|_{H^j(\mathbb R^n)}
\]
for every $v\in H^j(\mathbb R^n)$, where $[a,J_\varepsilon]:=aJ_\varepsilon-J_\varepsilon a$. The same assertion holds componentwise for matrix coefficients.
\end{lemma}

\begin{proof}
Begin with $j=0$ and fix $\ell\in\{1,\ldots,n\}$. For $v\in C_c^\infty(\mathbb R^n)$,
\[
[a,J_\varepsilon]v(x)
=
\int_{\mathbb R^n}
\bigl(a(x)-a(x-y)\bigr)\rho_\varepsilon(y)v(x-y)\,dy.
\]
The fundamental theorem of calculus for Banach-valued functions~\ref{teo:b5-fundamental-calculo-riemann-banach}, applied to $t\mapsto a(x-ty)$, gives
\[
|a(x)-a(x-y)|
\leq
\|Da\|_{L^\infty(\mathbb R^n)}\|y\|.
\]
By Young's inequality from Theorem~\ref{teo:young-convoluciones},
\[
\|[a,J_\varepsilon]v\|_{L^2(\mathbb R^n)}
\leq
\|Da\|_{L^\infty(\mathbb R^n)}
\bigl\|\|y\|\rho_\varepsilon(y)\bigr\|_{L^1(\mathbb R^n)}
\|v\|_{L^2(\mathbb R^n)}
\leq
C\varepsilon\|Da\|_{L^\infty(\mathbb R^n)}\|v\|_{L^2(\mathbb R^n)}.
\]
Differentiating the identity $[a,J_\varepsilon]v=aJ_\varepsilon v-J_\varepsilon(av)$ yields
\[
D_\ell[a,J_\varepsilon]v
=
(D_\ell a)J_\varepsilon v
+
\int_{\mathbb R^n}
\bigl(a(x)-a(x-y)\bigr)
(D_\ell\rho_\varepsilon)(y)v(x-y)\,dy.
\]
The first term is bounded in $L^2(\mathbb R^n)$ by $\|D_\ell a\|_{L^\infty(\mathbb R^n)}\|v\|_{L^2(\mathbb R^n)}$. For the second, use Theorem~\ref{teo:b5-fundamental-calculo-riemann-banach} again and observe that
\[
\int_{\mathbb R^n}\|y\|\,|D_\ell\rho_\varepsilon(y)|\,dy
=
\int_{\mathbb R^n}\|z\|\,|D_\ell\rho(z)|\,dz,
\]
which is independent of $\varepsilon$. This gives
\[
\|[a,J_\varepsilon]v\|_{H^1(\mathbb R^n)}
\leq
C_a\|v\|_{L^2(\mathbb R^n)}.
\]
Density extends the inequality to all of $v\in L^2(\mathbb R^n)$.

Proceed by induction. Suppose that for some $q\in\mathbb N_0$, every coefficient $b\in C_c^\infty(\mathbb R^n)$ satisfies
\[
\|[b,J_\varepsilon]w\|_{H^{q+1}(\mathbb R^n)}
\leq
C_{b,q}\|w\|_{H^q(\mathbb R^n)}.
\]
If $v\in H^{q+1}(\mathbb R^n)$, the distributional identity
\[
D_\ell[a,J_\varepsilon]v
=
[D_\ell a,J_\varepsilon]v
+[a,J_\varepsilon]D_\ell v
\]
and the induction hypothesis, applied to $D_\ell a$ and $a$, give
\[
\|D_\ell[a,J_\varepsilon]v\|_{H^{q+1}(\mathbb R^n)}
\leq
C\|v\|_{H^{q+1}(\mathbb R^n)}.
\]
Together with the bound in $L^2(\mathbb R^n)$, this yields
\[
\|[a,J_\varepsilon]v\|_{H^{q+2}(\mathbb R^n)}
\leq
C\|v\|_{H^{q+1}(\mathbb R^n)}.
\]
By induction, the assertion holds for every $j\in\mathbb N_0$.

Now let $j=-q$, with $q\in\mathbb N$. Since $J_\varepsilon$ is self-adjoint,
\[
[a,J_\varepsilon]^*
=-[\overline a,J_\varepsilon].
\]
The part already proved, with index $q-1$, states that $[\overline a,J_\varepsilon]\colon H^{q-1}(\mathbb R^n)\longrightarrow H^q(\mathbb R^n)$ is uniformly continuous. By duality, $[a,J_\varepsilon]\colon H^{-q}(\mathbb R^n)\longrightarrow H^{-q+1}(\mathbb R^n)$ is uniformly continuous. This completes the proof for every $j\in\mathbb Z$.
\end{proof}

\begin{corollary}[Commutator of a differential operator with a mollifier]
\label{cor:regularidad-intermedia-conmutador-pdo-molificador}
Let
\[
Q=\displaystyle\sum_{\substack{\alpha\in\mathbb N_0^n\\|\alpha|\leq m}}A_\alpha(x)D^\alpha
\]
be a differential system on $\mathbb R^n$ whose matrix coefficients belong to $C_c^\infty$. For every $j\in\mathbb Z$, there is a constant $C_{Q,j}>0$, independent of $0<\varepsilon\leq1$, such that
\[
\|[Q,J_\varepsilon]v\|_{H^{j-m+1}(\mathbb R^n)}
\leq
C_{Q,j}\|v\|_{H^j(\mathbb R^n)}
\]
for every $v\in H^j(\mathbb R^n,\mathbb K^r)$.
\end{corollary}

\begin{proof}
Constant-coefficient derivatives commute with convolution, so
\[
[Q,J_\varepsilon]v
=
\displaystyle\sum_{\substack{\alpha\in\mathbb N_0^n\\|\alpha|\leq m}}
[A_\alpha,J_\varepsilon]D^\alpha v.
\]
By Proposition~\ref{prop:regularidad-intermedia-propiedades-escala-hilbertiana}, $D^\alpha v\in H^{j-|\alpha|}(\mathbb R^n)$. Lemma~\ref{lem:regularidad-intermedia-conmutador-friedrichs} gives
\[
[A_\alpha,J_\varepsilon]D^\alpha v
\in
H^{j-|\alpha|+1}(\mathbb R^n).
\]
Since $|\alpha|\leq m$, the inclusion $H^{j-|\alpha|+1}(\mathbb R^n)\hookrightarrow H^{j-m+1}(\mathbb R^n)$ is continuous. Summing the finitely many terms gives the stated estimate.
\end{proof}

\section{Slobodeckij spaces}
\label{sec:espacios-slobodeckij}

Bessel potentials measure regularity through a Fourier multiplier. A second approach, intrinsic to the spatial variable, compares the values of a function at two points. The power in the denominator is chosen to make the expression invariant under Euclidean scaling. This construction leads to Sobolev--Slobodeckij spaces and remains meaningful on an open subset without first extending the function outside it.

\begin{definition}[Gagliardo seminorm and Slobodeckij space]
\label{def:seminorma-gagliardo-slobodeckij}
\glsadd{espacio-slobodeckij}
\index{Gagliardo seminorm}
\index{Slobodeckij space}
Let $\Omega\subseteq\mathbb R^n$ be open, and let $0<\sigma<1$ and $1\leq p<\infty$. For $u\in L^p(\Omega)$, define
\[
[u]_{W^{\sigma,p}(\Omega)}
:=
\left(
\int_\Omega\int_\Omega
\frac{|u(x)-u(y)|^p}{\|x-y\|^{n+\sigma p}}\,dx\,dy
\right)^{\frac{1}{p}}.
\]
The space $W^{\sigma,p}(\Omega)$ consists of functions for which this seminorm is finite and is equipped with
\[
\|u\|_{W^{\sigma,p}(\Omega)}
:=
\|u\|_{L^p(\Omega)}+[u]_{W^{\sigma,p}(\Omega)}.
\]
If $s=m+\sigma$, where $m\in\mathbb N_0$ and $0<\sigma<1$, set
\[
W^{s,p}(\Omega)
:=
\left\{
u\in W^{m,p}(\Omega)
\;\middle|\;
D^\alpha u\in W^{\sigma,p}(\Omega)
\text{ for }|\alpha|=m
\right\},
\]
with the norm
\[
\|u\|_{W^{s,p}(\Omega)}
:=
\|u\|_{W^{m,p}(\Omega)}
+
\displaystyle\sum_{\substack{\alpha\in\mathbb N_0^n\\|\alpha|=m}}[D^\alpha u]_{W^{\sigma,p}(\Omega)}.
\]
For $s=m\in\mathbb N_0$, retain the usual definition of $W^{m,p}(\Omega)$.
\end{definition}

The definition at integers is not obtained by formally setting $\sigma=0$ or $\sigma=1$: the corresponding integrals do not define the integer-order norm. We shall see later that this apparent discontinuity reflects a real difference between the fine indices of the Besov and Triebel--Lizorkin scales.

\begin{proposition}[Increment formulation]
\label{prop:gagliardo-mediante-incrementos}
If $u\in L^p(\mathbb R^n)$, $0<\sigma<1$, and $1\leq p<\infty$, then
\[
[u]_{W^{\sigma,p}(\mathbb R^n)}^p
=
\int_{\mathbb R^n}
\frac{\|\tau_hu-u\|_{L^p(\mathbb R^n)}^p}{\|h\|^{n+\sigma p}}\,dh,
\qquad
\tau_hu(x):=u(x+h).
\]
In particular, the functional on the right is translation invariant.
\end{proposition}

\begin{proof}
The integrand is nonnegative. Theorem~\ref{teo:tonelli} and the change of variables $y=x+h$ give
\[
\int_{\mathbb R^n}\int_{\mathbb R^n}
\frac{|u(x+h)-u(x)|^p}{\|h\|^{n+\sigma p}}\,dx\,dh
=
\int_{\mathbb R^n}\int_{\mathbb R^n}
\frac{|u(y)-u(x)|^p}{\|y-x\|^{n+\sigma p}}\,dy\,dx.
\]
Neither integral needs to be assumed finite beforehand.
\end{proof}

\begin{theorem}[Functional structure and density]
\label{teo:estructura-W-fraccionario}
Let $s>0$, $s\notin\mathbb N$, and $1\leq p<\infty$.
\begin{enumerate}[label=(\alph*)]
\item $W^{s,p}(\Omega)$ is a Banach space for every open subset $\Omega\subseteq\mathbb R^n$.
\item If $1<p<\infty$, then $W^{s,p}(\Omega)$ is reflexive; for every $1\leq p<\infty$, it is separable.
\item $C_c^\infty(\mathbb R^n)$ is dense in $W^{s,p}(\mathbb R^n)$.
\end{enumerate}
\end{theorem}

\begin{proof}
First treat $0<s<1$. Let
\[
d\mu_s(x,y):=\|x-y\|^{-n-sp}\,dx\,dy
\]
on $\Omega\times\Omega\setminus\{x=y\}$. The map
\[
\mathcal J\colon u\longmapsto
\left(u,\frac{u(x)-u(y)}{\|x-y\|^{\frac{n}{p}+s}}\right)
\]
identifies $W^{s,p}(\Omega)$ isometrically with a subspace of
\[
L^p(\Omega)\times L^p(\Omega\times\Omega,dx\,dy)
\]
when the product carries the sum norm. If instead the product carries the norm $\ell^p(\{1,2\})$, and we set $a:=\|u\|_{L^p(\Omega)}$ and $b:=[u]_{W^{s,p}(\Omega)}$, then
\[
 (a^p+b^p)^{\frac{1}{p}}\leq a+b
 \leq2^{1-\frac{1}{p}}(a^p+b^p)^{\frac{1}{p}}.
\]
This subspace is closed. Indeed, if $u_j\to u$ in $L^p(\Omega)$ and the quotients converge to $v$ in $L^p(\Omega\times\Omega)$, a subsequence converges pointwise almost everywhere in both spaces. Outside a null set, we obtain
\[
v(x,y)=\frac{u(x)-u(y)}{\|x-y\|^{\frac{n}{p}+s}}.
\]
Completeness of the product proves (a). The same argument, together with reflexivity and separability of the $L^p(\Omega)$ spaces established in Theorem~\ref{teo:dualidad-reflexividad-Lp}, proves (b). For $s=m+\sigma$, apply the argument to weak derivatives of order $m$ and use the fact that the operator $u\mapsto(D^\alpha u)_{|\alpha|\leq m}$ has closed graph in the product of $L^p(\Omega)$ spaces.

We prove (c). Choose $\eta\in C_c^\infty(\mathbb R^n)$ with $\eta=1$ on $B_{\mathrm{euc}}(0,1)$ and set $\eta_R(x)=\eta(\displaystyle\frac{x}{R})$. For $0<s<1$,
\[
\begin{split}
&|(1-\eta_R(x))u(x)-(1-\eta_R(y))u(y)|\\
&\qquad\leq
|1-\eta_R(x)|\,|u(x)-u(y)|
+|u(y)|\,|\eta_R(x)-\eta_R(y)|.
\end{split}
\]
The first term converges to zero in the weighted product measure by Theorem~\ref{convergencia dominada}. For the second, use
\[
|\eta_R(x)-\eta_R(y)|
\leq
C\min\left\{1,\frac{\|x-y\|}{R}\right\}.
\]
After the substitution $h=x-y$,
\[
\int_{\mathbb R^n}
\frac{\min\left\{1,\frac{\|h\|}{R}\right\}^p}{\|h\|^{n+sp}}\,dh
=
C_{n,s,p}R^{-sp};
\]
so the second term tends to zero. We conclude that $\eta_Ru\to u$ in $W^{s,p}(\mathbb R^n)$.

Now let $\rho$ be a standard mollifier and $u_\varepsilon=\rho_\varepsilon*u$. Theorem~\ref{teo:minkowski-integral} and translation invariance give
\[
\|u_\varepsilon-u\|_{W^{s,p}(\mathbb R^n)}
\leq
\int_{\mathbb R^n}\rho(z)\|\tau_{\varepsilon z}u-u\|_{W^{s,p}(\mathbb R^n)}\,dz.
\]
Translations are continuous in $W^{s,p}(\mathbb R^n)$: this follows by applying Proposition~\ref{prop:continuidad-traslaciones-Lp} both to $u\in L^p(\mathbb R^n)$ and to the two-variable function representing its difference quotient in $L^p(\mathbb R^{2n})$. Theorem~\ref{convergencia dominada} implies $u_\varepsilon\to u$ in $W^{s,p}(\mathbb R^n)$. Apply the cutoff first, then regularization; the result belongs to $C_c^\infty(\mathbb R^n)$. For \(s=m+\sigma\), \(0<\sigma<1\), also expand the terms produced by the cutoff. If \(|\alpha|=m\), the Leibniz rule gives
\[
D^\alpha(\eta_Ru-u)
=(\eta_R-1)D^\alpha u+
\sum_{\substack{\beta\in\mathbb N_0^n\\
\beta\leq\alpha,\ \beta\neq\alpha}}
\binom{\alpha}{\beta}D^{\alpha-\beta}\eta_R\,D^\beta u.
\]
The first term converges to zero in \(W^{\sigma,p}(\mathbb R^n)\) by what has already been proved. For the remaining terms, \(D^\beta u\in W^{1,p}(\mathbb R^n)\). The increment lemma, whose proof uses only integer-order Sobolev properties, gives
\[
\|w\|_{W^{\sigma,p}(\mathbb R^n)}\leq C\|w\|_{W^{1,p}(\mathbb R^n)},
\qquad w\in W^{1,p}(\mathbb R^n):
\]
in the increment integral, use
\(\|\tau_hw-w\|_{L^p(\mathbb R^n)}\leq\|h\|\|\nabla w\|_{L^p(\mathbb R^n)}\)
for \(\|h\|<1\), and \(2\|w\|_{L^p(\mathbb R^n)}\) on the complement. If \(\gamma=\alpha-\beta\neq0\), both \(\|D^\gamma\eta_R\|_{L^\infty(\mathbb R^n)}\) and \(\|D(D^\gamma\eta_R)\|_{L^\infty(\mathbb R^n)}\) tend to zero. The formula for the difference of a product and these two bounds imply
\[
\|D^\gamma\eta_R\,D^\beta u\|_{W^{\sigma,p}(\mathbb R^n)}
\leq C\bigl(\|D^\gamma\eta_R\|_{L^\infty(\mathbb R^n)}
+\|D(D^\gamma\eta_R)\|_{L^\infty(\mathbb R^n)}\bigr)
\|D^\beta u\|_{W^{\sigma,p}(\mathbb R^n)}\longrightarrow0.
\]
Lower-order derivatives are controlled by the same Leibniz rule in \(L^p(\mathbb R^n)\). This justifies the cutoff in the full norm. Regularization commutes with each weak derivative, and the continuity of translations already proved applies to the quotients of derivatives of order \(m\). Thus regularization after the cutoff establishes density also at higher orders.
\end{proof}

\subsection{The Slobodeckij interpolation identity}

The Gagliardo seminorm is not an isolated choice: it is exactly the quantity measuring the position of $u$ between $L^p(\mathbb R^n)$ and $W^{1,p}(\mathbb R^n)$. This interpretation explains both the exponent $n+sp$ and many subsequent properties.

\begin{lemma}[Increments of a Sobolev function]
\label{lem:incrementos-W1p}
If $w\in W^{1,p}(\mathbb R^n)$ and $1\leq p<\infty$, then
\[
\|\tau_hw-w\|_{L^p(\mathbb R^n)}
\leq
\|h\|\,\|\nabla w\|_{L^p(\mathbb R^n)}
\qquad\text{for every }h\in\mathbb R^n.
\]
\end{lemma}

\begin{proof}
For $w\in C_c^\infty(\mathbb R^n)$, Theorem~\ref{teo:b5-fundamental-calculo-riemann-banach} gives
\[
w(x+h)-w(x)
=
\int_0^1\nabla w(x+th)\cdot h\,dt.
\]
Theorem~\ref{teo:minkowski-integral} and translation invariance of the $L^p(\mathbb R^n)$ norm yield the inequality. Density of $C_c^\infty(\mathbb R^n)$ in $W^{1,p}(\mathbb R^n)$, Proposition~\ref{Cc es denso en W1p}, allows passage to every $w\in W^{1,p}(\mathbb R^n)$.
\end{proof}

\begin{theorem}[Slobodeckij as a real interpolation space]
\label{teo:slobodeckij-interpolacion-real}
\index{real interpolation!Slobodeckij spaces}
Let $0<s<1$ and $1\leq p<\infty$. Then
\[
\bigl(L^p(\mathbb R^n),W^{1,p}(\mathbb R^n)\bigr)_{s,p}
=
W^{s,p}(\mathbb R^n)
\]
and there are constants $c_{n,s,p},C_{n,s,p}>0$ such that
\[
c_{n,s,p}\|u\|_{W^{s,p}(\mathbb R^n)}
\leq
\|u\|_{(L^p(\mathbb R^n),W^{1,p}(\mathbb R^n))_{s,p}}
\leq
C_{n,s,p}\|u\|_{W^{s,p}(\mathbb R^n)}
\]
for every $u\in W^{s,p}(\mathbb R^n)$.
\end{theorem}

\begin{proof}
Write $K(t,u)$ for the functional of the pair $(L^p(\mathbb R^n),W^{1,p}(\mathbb R^n))$. If $u=v+w$, with $v\in L^p(\mathbb R^n)$ and $w\in W^{1,p}(\mathbb R^n)$, Lemma~\ref{lem:incrementos-W1p} and the elementary inequality $(a+b)^p\leq2^{p-1}(a^p+b^p)$ give
\[
\|\tau_hu-u\|_{L^p(\mathbb R^n)}^p
\leq
C_p\bigl(\|v\|_{L^p(\mathbb R^n)}+\|h\|\|w\|_{W^{1,p}(\mathbb R^n)}\bigr)^p.
\]
Taking the infimum and using polar coordinates,
\[
[u]_{W^{s,p}(\mathbb R^n)}^p
\leq
C\int_{\mathbb R^n}
\frac{K(\|h\|,u)^p}{\|h\|^{n+sp}}\,dh
=
C|\mathbb S^{n-1}|\int_0^\infty
\bigl(t^{-s}K(t,u)\bigr)^p\,\frac{dt}{t}.
\]
The embedding of the interpolated space into $L^p(\mathbb R^n)+W^{1,p}(\mathbb R^n)=L^p(\mathbb R^n)$, proved in Theorem~\ref{teo:propiedades-basicas-metodo-K}, also controls the $L^p(\mathbb R^n)$ norm. This proves one inclusion.

For the converse, take a nonnegative mollifier $\rho$ supported in the unit ball and set $u_t=\rho_t*u$. The decomposition
\[
u=(u-u_t)+u_t
\]
is admissible for $K(t,u)$. Proposition~\ref{desigualdad de holder}, with respect to the probability measure $\rho_t(h)\,dh$, gives
\[
\|u-u_t\|_{L^p(\mathbb R^n)}^p
\leq
Ct^{-n}\int_{\{h\in\mathbb R^n\mid\,\|h\|<t\}}\|\tau_hu-u\|_{L^p(\mathbb R^n)}^p\,dh.
\]
Since $\displaystyle\int_{\mathbb R^n}D_i\rho_t(h)\,dh=0$,
\[
D_iu_t(x)
=
\int_{\mathbb R^n}D_i\rho_t(h)\bigl(u(x-h)-u(x)\bigr)\,dh,
\]
and another application of Proposition~\ref{desigualdad de holder} gives
\[
\|\nabla u_t\|_{L^p(\mathbb R^n)}^p
\leq
Ct^{-n-p}\int_{\{h\in\mathbb R^n\mid\,\|h\|<t\}}\|\tau_hu-u\|_{L^p(\mathbb R^n)}^p\,dh.
\]
Moreover, Theorem~\ref{teo:young-convoluciones} gives $\|u_t\|_{L^p(\mathbb R^n)}\leq\|u\|_{L^p(\mathbb R^n)}$. Consequently,
\[
K(t,u)^p
\leq
C\left[
t^{-n}\int_{\{\|h\|<t\}}\|\tau_hu-u\|_{L^p(\mathbb R^n)}^p\,dh
+t^p\|u\|_{L^p(\mathbb R^n)}^p
\right]
\]
for $0<t<1$. Multiply by $t^{-sp}$, integrate with respect to $\displaystyle\frac{dt}{t}$, and apply Theorem~\ref{teo:tonelli}. The first term gives
\[
\int_{\{\|h\|<1\}}\|\tau_hu-u\|_{L^p(\mathbb R^n)}^p
\left(
\int_{\|h\|}^1t^{-n-sp-1}\,dt
\right)dh
\leq
C[u]_{W^{s,p}(\mathbb R^n)}^p.
\]
The second term is $C\|u\|_{L^p(\mathbb R^n)}^p\displaystyle\int_0^1t^{(1-s)p}\,\displaystyle\frac{dt}{t}$. Apply the two quantified inequalities for the truncated norm in Theorem~\ref{teo:propiedades-basicas-metodo-K} and absorb their constants, which depend only on $s$ and $p$; this gives
\[
\|u\|_{(L^p(\mathbb R^n),W^{1,p}(\mathbb R^n))_{s,p}}
\leq
C\|u\|_{W^{s,p}(\mathbb R^n)}.
\]
\end{proof}

\begin{corollary}[Interpolation between consecutive integer orders]
\label{cor:slobodeckij-interpolacion-altos-ordenes}
Let $m\in\mathbb N_0$, $0<s<1$, and $1<p<\infty$. Then
\[
\bigl(W^{m,p}(\mathbb R^n),W^{m+1,p}(\mathbb R^n)\bigr)_{s,p}
=
W^{m+s,p}(\mathbb R^n)
\]
as sets, and there are constants $c_{n,m,s,p},C_{n,m,s,p}>0$ such that
\[
 c_{n,m,s,p}\|u\|_{W^{m+s,p}(\mathbb R^n)}
 \leq
 \|u\|_{(W^{m,p}(\mathbb R^n),W^{m+1,p}(\mathbb R^n))_{s,p}}
 \leq
 C_{n,m,s,p}\|u\|_{W^{m+s,p}(\mathbb R^n)}
\]
for every $u\in W^{m+s,p}(\mathbb R^n)$.
\end{corollary}

\begin{proof}
For \(m=0\), this is the preceding theorem. Assume \(m\geq1\). If \(u=v+w\), with \(v\in W^{m,p}(\mathbb R^n)\) and \(w\in W^{m+1,p}(\mathbb R^n)\), Lemma~\ref{lem:incrementos-W1p} gives, for \(\alpha\in\mathbb N_0^n\), \(|\alpha|=m\),
\[
\|\tau_hD^\alpha u-D^\alpha u\|_{L^p(\mathbb R^n)}
\leq 2\|v\|_{W^{m,p}(\mathbb R^n)}+
C\|h\|\|w\|_{W^{m+1,p}(\mathbb R^n)}.
\]
Take the infimum over these decompositions and integrate with weight \(\|h\|^{-n-sp}\) exactly as in the proof for \(m=0\). This controls all seminorms of order \(m+s\). Moreover, \(W^{m+1,p}(\mathbb R^n)\hookrightarrow W^{m,p}(\mathbb R^n)\) implies
\[
K(1,u;W^{m,p}(\mathbb R^n),W^{m+1,p}(\mathbb R^n))\geq c\|u\|_{W^{m,p}(\mathbb R^n)}.
\]
Monotonicity of \(K\) and the portion of the interpolation integral corresponding to \(1\leq t\leq2\) also control this norm. We obtain the first continuous inclusion.

For the converse, fix the mollifier \(\rho\) from the preceding proof, supported in the unit ball, and write \(u=(u-u_t)+u_t\), where \(u_t=\rho_t*u\). If \(|\alpha|<m\), the increment lemma gives
\[
\|D^\alpha(u-u_t)\|_{L^p(\mathbb R^n)}
\leq Ct\|u\|_{W^{m,p}(\mathbb R^n)},\qquad 0<t<1.
\]
For \(|\alpha|=m\), Minkowski's and Hölder's inequalities give
\[
\|D^\alpha(u-u_t)\|_{L^p(\mathbb R^n)}^p
\leq Ct^{-n}\int_{\|h\|<t}
\|\tau_hD^\alpha u-D^\alpha u\|_{L^p(\mathbb R^n)}^p\,dh.
\]
Since \(\displaystyle \int_{\mathbb R^n}D_\ell\rho_t=0\), also
\[
t^p\|D_\ell D^\alpha u_t\|_{L^p(\mathbb R^n)}^p
\leq Ct^{-n}\int_{\|h\|<t}
\|\tau_hD^\alpha u-D^\alpha u\|_{L^p(\mathbb R^n)}^p\,dh,
\qquad \ell\in\{1,\ldots,n\}.
\]
Derivatives of \(u_t\) of order at most \(m\) are controlled by Young's inequality, so
\[
\begin{split}
K(t,u;W^{m,p}(\mathbb R^n),W^{m+1,p}(\mathbb R^n))^p
\leq C t^p\|u\|_{W^{m,p}(\mathbb R^n)}^p
+C t^{-n}\int_{\|h\|<t}
\sum_{\substack{\alpha\in\mathbb N_0^n\\|\alpha|=m}}
\|\tau_hD^\alpha u-D^\alpha u\|_{L^p(\mathbb R^n)}^p\,dh.
\end{split}
\]
Multiply by \(t^{-sp}\) and integrate with \(\displaystyle\frac{dt}{t}\) over \(0<t<1\). The first term is integrable because \(s<1\). For the second, Tonelli's theorem and
\[
\int_{\|h\|}^1t^{-n-sp-1}\,dt
\leq \frac{\|h\|^{-n-sp}}{n+sp},
\qquad 0<\|h\|<1,
\]
yield the Gagliardo seminorms of derivatives of order \(m\). For \(t\geq1\), use the decomposition \(u=u+0\); the remaining integral is bounded by \(\displaystyle\frac{\|u\|_{W^{m,p}(\mathbb R^n)}^p}{sp}\). This proves the second inclusion.

For $r\in\{m,m+1\}$, the integer-order norms used are the \(\ell^p(\{\alpha\in\mathbb N_0^n\mid|\alpha|\leq r\})\) norms of this book. If a sum of norms appears in an estimate, the comparison is explicit: for \(r\in\{m,m+1\}\), with \(N_r=\#\{\alpha\in\mathbb N_0^n\mid|\alpha|\leq r\}\),
\[
\left(\sum_{\substack{\alpha\in\mathbb N_0^n\\|\alpha|\leq r}}
a_\alpha^p\right)^{\frac{1}{p}}
\leq
\sum_{\substack{\alpha\in\mathbb N_0^n\\|\alpha|\leq r}}a_\alpha
\leq N_r^{1-\frac{1}{p}}
\left(\sum_{\substack{\alpha\in\mathbb N_0^n\\|\alpha|\leq r}}
a_\alpha^p\right)^{\frac{1}{p}}.
\]
The constants therefore depend only on \(n,m,s,p\) and the fixed mollifier. The same argument holds for \(p=1\).
\end{proof}

\section{Besov and Triebel--Lizorkin spaces}
\label{sec:besov-triebel-lizorkin-euclideos}

The Littlewood--Paley decomposition from Chapter~\ref{cap:distribuciones-fourier} allows a distribution to be studied by separating its content into frequency scales. Each block $\Delta_j u$ measures the part of $u$ lying roughly at frequencies of size $2^j$, while the block $\Delta_0u$ collects low frequencies. Once this separation is made, there are two natural ways to combine the information. In the first, one takes the spatial norm of each block and then sums over scales; in the second, one first combines all blocks pointwise and only then takes the spatial norm. The first procedure leads to Besov spaces, the second to Triebel--Lizorkin spaces.

The difference between these two orders of combination is not merely notational. When the indices agree, Tonelli's theorem permits interchange of sum and integral; otherwise, the position of the sequence norm generally changes the resulting space. This distinction will later explain why Sobolev--Slobodeckij spaces are identified with $B^s_{p,p}(\mathbb R^n)$, while Bessel potentials are identified with $F^s_{p,2}(\mathbb R^n)$. This organization of the scales is the classical Littlewood--Paley description; see, for example, \cite[Chapters~1--2]{Triebel2} and \cite[Chapter~2]{Sawano2018}.

We use the resolution $(\phi_j)_{j\in\mathbb N_0}$ and the blocks $\Delta_j$ from Proposition~\ref{prop:resolucion-diadica-unidad} and Definition~\ref{def:bloques-diadicos}. Also choose real radial functions $\widetilde\phi_j\in C_c^\infty(\mathbb R^n)$, $j\in\mathbb N_0$, equal to one on a neighborhood of $\operatorname{supp}(\phi_j)$ and having uniform dyadic supports. Write
\[
\widetilde\Delta_j:=\mathcal F^{-1}\widetilde\phi_j\mathcal F,
\qquad j\in\mathbb N_0.
\]
Then $\widetilde\Delta_j\Delta_j=\Delta_j$ for every $j\in\mathbb N_0$.

\subsection{Multiplier tools}

The dyadic construction uses two estimates: continuity of certain Fourier multipliers and control of a function by its frequency blocks. Their constants will need to be uniform in the block index.

\begin{theorem}[Mikhlin multipliers]
\label{teo:multiplicadores-mikhlin}
\index{Mikhlin theorem}
Let $1<p<\infty$ and let $N\geq n+2$ be an integer. Suppose that $m\in C^N(\mathbb R^n\setminus\{0\})$ satisfies
\[
M_N(m)
:=
\max_{\substack{\alpha\in\mathbb N_0^n\\|\alpha|\leq N}}
\sup_{\xi\in\mathbb R^n\setminus\{0\}}\|\xi\|^{|\alpha|}|D^\alpha m(\xi)|
<\infty.
\]
Then the operator initially defined on $\mathcal S(\mathbb R^n)$ by
\[
T_mf:=\mathcal F^{-1}(m\widehat f)
\]
extends uniquely to a continuous operator on $L^p(\mathbb R^n)$, and
\[
\|T_mf\|_{L^p(\mathbb R^n)}
\leq
C_{n,p,N}M_N(m)\|f\|_{L^p(\mathbb R^n)}.
\]
The same bound holds for a family of symbols with uniform Mikhlin constants.
\end{theorem}

\begin{proof}
If \(M_N(m)=0\), the symbol vanishes. Otherwise, first divide \(m\) by \(M_N(m)\) and restore this factor at the end. We may thus assume \(M_N(m)=1\). Let $\vartheta\in C_c^\infty(\mathbb R^n\setminus\{0\})$ satisfy $\displaystyle\sum_{j\in\mathbb Z}\vartheta(2^{-j}\xi)=1$ for $\xi\neq0$, and set $m_j(\xi)=m(\xi)\vartheta(2^{-j}\xi)$. The substitution $\xi=2^j\zeta$ and the hypotheses show that
\[
\max_{\substack{\alpha\in\mathbb N_0^n\\|\alpha|\leq N}}
\|D^\alpha[m(2^j\,\cdot)\vartheta]\|_{L^\infty(\mathbb R^n)}
\leq
CM_N(m)
\]
uniformly in $j$. Integration by parts in the inverse transform gives, for the kernel $K_j=\mathcal F^{-1}m_j$,
\[
|K_j(x)|
\leq
CM_N(m)2^{jn}(1+2^j\|x\|)^{-N}
\]
and applying the same computation to a derivative gives
\[
|\nabla K_j(x)|
\leq
CM_N(m)2^{j(n+1)}(1+2^j\|x\|)^{-N}.
\]
The sum $K=\displaystyle\sum_{j\in\mathbb Z}K_j$ converges locally uniformly on $\mathbb R^n\setminus\{0\}$. Separate indices according to whether $2^j\|x\|\leq1$ or $2^j\|x\|>1$, and sum the two geometric series; this gives
\[
|K(x)|\leq C M_N(m)\|x\|^{-n},
\qquad
|K(x-h)-K(x)|
\leq
C M_N(m)\frac{\|h\|}{\|x\|^{n+1}}
\]
when $\|h\|\leq\displaystyle\frac{1}{2}\|x\|$. The operator is bounded on $L^2(\mathbb R^n)$ by Theorem~\ref{teo:plancherel-fourier}, since $m$ is bounded. With the unitary Fourier normalization, the convolution kernel is \((2\pi)^{-\displaystyle\frac{n}{2}}K\). To make this representation precise, let $b\in L^1\cap L^2$
have compact support and set $m^{(L)}:=\sum_{j=-L}^{L}m_j$.
The partial symbols are uniformly bounded by finite overlap of
the annuli and converge to $m$ almost everywhere. Dominated
convergence in frequency and Plancherel's theorem give
\[
T_{m^{(L)}}b\longrightarrow T_mb\quad\text{in }L^2.
\]
On each compact set disjoint from $\operatorname{supp}b$, the
kernel series converges absolutely and uniformly on the differences
$x-y$, with $y\in\operatorname{supp}b$.
Hence the integrals of its partial sums converge to
$(2\pi)^{-\displaystyle\frac{n}{2}}\int K(x-y)b(y)\,dy$.
A subsequence of the $L^2$ convergence converges almost everywhere;
comparing these limits identifies this integral with $T_mb$
almost everywhere outside the support.

Let \(f\in L^1(\mathbb R^n)\cap L^2(\mathbb R^n)\), let \(\lambda>0\), and apply the Calderón--Zygmund decomposition from Theorem~\ref{teo:descomposicion-calderon-zygmund}:
\[
f=g+\sum_{Q\in\mathcal Q}b_Q,\qquad
\|g\|_{L^2(\mathbb R^n)}^2\leq C\lambda\|f\|_{L^1(\mathbb R^n)},\qquad
\sum_{Q\in\mathcal Q}\|b_Q\|_{L^1(\mathbb R^n)}\leq C\|f\|_{L^1(\mathbb R^n)}.
\]
Each \(b_Q\) has integral zero, vanishes outside \(Q\),
and its topological support is contained in \(\overline Q\). Moreover, \(\displaystyle \sum_{Q\in\mathcal Q}|Q|\leq C\displaystyle\frac{\|f\|_{L^1(\mathbb R^n)}}{\lambda}\). The decomposition also converges in $L^2$ for this $f$.
Indeed, $\|g\|_\infty\leq C\lambda$ and
$\|g\|_1\leq C\|f\|_1$ imply $g\in L^2$ and the stated bound.
Thus $b:=f-g\in L^2$. The half-open cubes are disjoint and
$b_Q=b\mathbf1_Q$, so Tonelli's theorem gives
\[
\sum_{Q\in\mathcal Q}\|b_Q\|_2^2=\|b\|_2^2<\infty.
\]
After enumerating the cubes, the squared $L^2$ norm of each tail
is the tail of this numerical series. Hence $\sum_Q b_Q=b$
in $L^2$, and continuity of $T_m$ on $L^2$ allows the operator
to be applied to this sum.

Write \(Q^*=4\sqrt n\,Q\), with the same center \(c_Q\). If \(x\notin Q^*\) and \(y\in Q\), then \(\|y-c_Q\|\leq\displaystyle\frac{\|x-c_Q\|}{2}\). By cancellation,
\[
T_mb_Q(x)=(2\pi)^{-\frac{n}{2}}\int_Q
[K(x-y)-K(x-c_Q)]b_Q(y)\,dy,
\]
and the kernel estimate implies
\[
|T_mb_Q(x)|
\leq C\frac{\ell(Q)}{\|x-c_Q\|^{n+1}}\|b_Q\|_{L^1(\mathbb R^n)}.
\]
The integral of \(\ell(Q)\|x-c_Q\|^{-n-1}\) outside \(Q^*\) is bounded by a dimensional constant. Thus, if \(\displaystyle E=\bigcup_{Q\in\mathcal Q}Q^*\),
\[
|E|\leq C\frac{\|f\|_{L^1(\mathbb R^n)}}{\lambda},\qquad
\int_{\mathbb R^n\setminus E}
\sum_{Q\in\mathcal Q}|T_mb_Q|
\leq C\|f\|_{L^1(\mathbb R^n)}.
\]
The partial sums converge in $L^2$ to $T_mb$. Outside $E$,
the preceding bound gives absolute convergence of the series
almost everywhere, and even convergence in $L^1$ on that region.
A subsequence of the partial sums converges almost everywhere
to its $L^2$ limit; uniqueness of the pointwise limit identifies
the series outside $E$ with $T_m(f-g)$.
Chebyshev's inequality and the \(L^2(\mathbb R^n)\) bound for the good part give
\[
\begin{split}
|\{x\in\mathbb R^n\mid |T_mf(x)|>\lambda\}|
&\leq |E|+
\frac{4}{\lambda^2}\|T_mg\|_{L^2(\mathbb R^n)}^2
+\frac{2}{\lambda}
\int_{\mathbb R^n\setminus E}|T_m(f-g)|\\
&\leq C\frac{\|f\|_{L^1(\mathbb R^n)}}{\lambda}.
\end{split}
\]
Theorem~\ref{teo:interpolacion-marcinkiewicz} interpolates the
weak $(1,1)$ bound with the strong $(2,2)$ bound and gives a
continuous operator on $L^p$ for $1<p<2$. It agrees with the
Plancherel operator on $\mathcal S$: for $f\in\mathcal S$, choose
simple functions of finite-measure support converging to $f$
simultaneously in $L^p$ and $L^2$. These are obtained by truncating
the support and values and approximating the values on sets of
finite measure. Their images converge in $L^p$ to the extension
and in $L^2$ to the original operator; the two limits agree almost
everywhere by uniqueness of limits in measure on balls.
The case $p=2$ is already covered by Plancherel's theorem.

If $p>2$, the symbol $\overline m$ satisfies the same bounds
and its operator is bounded on $L^{p'}$. For $f,h\in\mathcal S$,
Plancherel's theorem and Hölder's inequality imply
\[
\left|\int T_mf\,\overline h\right|
=\left|\int f\,\overline{T_{\overline m}h}\right|
\leq C\|f\|_p\|h\|_{p'}.
\]
The dual criterion in Lemma~\ref{lem:criterio-dual-Lp-interpolacion}
gives $T_mf\in L^p$ with the same bound; it applies because
$T_mf\in L^2\subseteq L^1_{\mathrm{loc}}$.
Density of $\mathcal S$ in $L^p$ yields the unique continuous
extension for all the exponents in the statement. Finally,
restore the factor $M_N(m)$ to obtain the asserted uniformity.
\end{proof}

\begin{theorem}[Littlewood--Paley square function]
\label{teo:funcion-cuadrada-littlewood-paley}
\index{Littlewood--Paley square function}
Let $1<p<\infty$. There are constants $c_p,C_p>0$ such that for every $f\in L^p(\mathbb R^n)$,
\[
c_p\|f\|_{L^p(\mathbb R^n)}
\leq
\left\|
\left(\sum_{j\in\mathbb N_0}|\Delta_jf|^2\right)^{\frac{1}{2}}
\right\|_{L^p(\mathbb R^n)}
\leq
C_p\|f\|_{L^p(\mathbb R^n)}.
\]
\end{theorem}

\begin{proof}
We first bound the square function by the norm of $f$. Let $J\subseteq\mathbb N_0$ be finite. For each choice of signs $(\varepsilon_j)_{j\in J}\in\{-1,1\}^J$, the symbol
\[
m_\varepsilon(\xi):=\sum_{j\in J}\varepsilon_j\phi_j(\xi)
\]
satisfies the hypotheses of Theorem~\ref{teo:multiplicadores-mikhlin} with constants independent of $J$ and the signs: finite interaction of the supports ensures that at each $\xi\neq0$ only a uniformly bounded number of terms occur, and after scaling, the derivatives of $\phi_j$ have the required bounds. Therefore,
\[
\left\|\sum_{j\in J}\varepsilon_j\Delta_jf\right\|_{L^p(\mathbb R^n)}
\leq
C_p\|f\|_{L^p(\mathbb R^n)}.
\]
Take $\varepsilon_j=r_j(t)$, where $(r_j)_{j\in\mathbb N_0}$ are the Rademacher functions, integrate in $t\in[0,1]$, and use Khintchine's theorem~\ref{teo:khintchine} pointwise in $x$. Theorem~\ref{teo:tonelli} allows the integrals to be interchanged and gives
\[
\left\|
\left(\sum_{j\in J}|\Delta_jf|^2\right)^{\frac{1}{2}}
\right\|_{L^p(\mathbb R^n)}
\leq
C_p\|f\|_{L^p(\mathbb R^n)}.
\]
Letting $J$ increase to $\mathbb N_0$ and applying monotone convergence yields the full upper bound.

For the reverse inequality, first take $f,g\in\mathcal S(\mathbb R^n)$. Since $\widetilde\Delta_j\Delta_j=\Delta_j$ and the dyadic resolution reconstructs the identity in $\mathcal S'(\mathbb R^n)$,
\[
\langle f,g\rangle
=
\sum_{j\in\mathbb N_0}\langle\Delta_jf,\widetilde\Delta_jg\rangle.
\]
Finite interaction of the supports allows the enlarged blocks to be chosen so that the series converges absolutely for Schwartz functions. Cauchy--Schwarz in $\ell^2(\mathbb N_0)$ gives
\[
|\langle f,g\rangle|
\leq
\int_{\mathbb R^n}
\left(\sum_{j\in\mathbb N_0}|\Delta_jf(x)|^2\right)^{\frac{1}{2}}
\left(\sum_{j\in\mathbb N_0}|\widetilde\Delta_jg(x)|^2\right)^{\frac{1}{2}}dx.
\]
Apply Hölder's inequality with exponents $p$ and $p'$ and the square-function bound already proved, now for the enlarged resolution and $g\in L^{p'}(\mathbb R^n)$. Thus,
\[
|\langle f,g\rangle|
\leq
C
\left\|
\left(\sum_{j\in\mathbb N_0}|\Delta_jf|^2\right)^{\frac{1}{2}}
\right\|_{L^p(\mathbb R^n)}
\|g\|_{L^{p'}(\mathbb R^n)}.
\]
Duality from Theorem~\ref{teo:dualidad-reflexividad-Lp} gives the reverse inequality for $f\in\mathcal S(\mathbb R^n)$. Finally, density of $\mathcal S(\mathbb R^n)$ in $L^p(\mathbb R^n)$, Proposition~\ref{prop:densidad-D-en-S-y-S-en-Lp}, and continuity of the square function extend the equivalence to all of $L^p(\mathbb R^n)$.
\end{proof}

For $a>0$, we shall use the abbreviation
\[
\mathcal M_af:=\bigl[\mathcal M(|f|^a)\bigr]^{\frac{1}{a}}.
\]

\begin{lemma}[Dyadic multipliers on sequences]
\label{lem:multiplicadores-diadicos-vectoriales}
Let $0<p<\infty$, $0<q\leq\infty$, and $\displaystyle 0<a<\min\{p,q\}$, with the convention that $q$ is omitted from this minimum if $q=\infty$. Suppose that $(f_j)_{j\in\mathbb N_0}\subseteq\mathcal S'(\mathbb R^n)$ is a family whose elements are represented by functions and satisfy, for every $j\in\mathbb N_0$,
\[
\operatorname{supp}(\widehat f_j)
\subseteq
\{\xi\in\mathbb R^n\mid \|\xi\|\leq A2^j\}.
\]
If $T_jf_j=K_j*f_j$, where $K_j(x)=2^{jn}K^{(j)}(2^jx)$ and
\[
\sup_{j\in\mathbb N_0}\sup_{x\in\mathbb R^n}
(1+\|x\|)^{N+n+1}|K^{(j)}(x)|<\infty
\]
for some $N>\displaystyle\frac{n}{a}$, then
\[
\left\|
\left(\displaystyle\sum_{j\in\mathbb N_0}|T_jf_j|^q\right)^{\frac{1}{q}}
\right\|_{L^p(\mathbb R^n)}
\leq
C
\left\|
\left(\displaystyle\sum_{j\in\mathbb N_0}|f_j|^q\right)^{\frac{1}{q}}
\right\|_{L^p(\mathbb R^n)},
\]
with the usual suprema when $q=\infty$.
\end{lemma}

\begin{proof}
Lemma~\ref{lem:control-maximal-peetre} gives $|T_jf_j|\leq C\mathcal M_af_j$ pointwise. If $q<\infty$, apply Theorem~\ref{teo:fefferman-stein-vectorial} to the sequence $(|f_j|^a)_{j\in\mathbb N_0}$ with exponents $\displaystyle\frac{p}{a}$ and $\displaystyle\frac{q}{a}$, both greater than one, and raise the result to the power $\displaystyle\frac{1}{a}$. If $q=\infty$, use
\[
\sup_{j\in\mathbb N_0}\mathcal M(|f_j|^a)
\leq
\mathcal M\left(\sup_{j\in\mathbb N_0}|f_j|^a\right)
\]
and continuity of $\mathcal M$ on $L^{\displaystyle\frac{p}{a}}(\mathbb R^n)$ from Theorem~\ref{teo:maximal-hardy-littlewood}.
\end{proof}

\begin{definition}[Besov and Triebel--Lizorkin spaces]
\label{def:besov-triebel-lizorkin}
\index{Besov space}
\index{Triebel--Lizorkin space}
Let $s\in\mathbb R$ and $1\leq q\leq\infty$. In this section, we work in the Banach range of these scales; classical theory extends the definitions to more general quasi-Banach indices, but these endpoints will not be needed for later results.
\begin{enumerate}[label=(\alph*)]
\item If $1\leq p\leq\infty$, define
\[
B^s_{p,q}(\mathbb R^n)
:=
\left\{
u\in\mathcal S'(\mathbb R^n)
\;\middle|\;
\|u\|_{B^s_{p,q}(\mathbb R^n)}<\infty
\right\},
\]
where
\[
\|u\|_{B^s_{p,q}(\mathbb R^n)}
:=
\left\|
\bigl(2^{js}\|\Delta_ju\|_{L^p(\mathbb R^n)}\bigr)_{j\in\mathbb N_0}
\right\|_{\ell^q(\mathbb N_0)}.
\]
\item If $1\leq p<\infty$, define
\[
F^s_{p,q}(\mathbb R^n)
:=
\left\{
u\in\mathcal S'(\mathbb R^n)
\;\middle|\;
\|u\|_{F^s_{p,q}(\mathbb R^n)}<\infty
\right\},
\]
where
\[
\|u\|_{F^s_{p,q}(\mathbb R^n)}
:=
\left\|
\left(\displaystyle\sum_{j\in\mathbb N_0}2^{jsq}|\Delta_ju|^q\right)^{\frac{1}{q}}
\right\|_{L^p(\mathbb R^n)},
\]
and, if \(q=\infty\),
\[
\|u\|_{F^s_{p,\infty}(\mathbb R^n)}
:=
\left\|\sup_{j\in\mathbb N_0}2^{js}|\Delta_ju|\right\|_{L^p(\mathbb R^n)}.
\]
\end{enumerate}
\end{definition}

\begin{definition}[Retraction and coretraction]
\label{def:retraccion-corretraccion-espacios-banach}
Let $X$ and $Y$ be Banach spaces. A continuous linear operator $R\colon X\longrightarrow Y$ is a \textit{retraction} if there is a continuous linear operator $S\colon Y\longrightarrow X$ such that
\[
RS=I_Y.
\]
In that case, $S$ is called a \textit{coretraction} of $R$ or a continuous right inverse. In particular, every retraction is surjective and every coretraction is injective.
\end{definition}

We also adopt the convention \(F^s_{\infty,\infty}(\mathbb R^n):=B^s_{\infty,\infty}(\mathbb R^n)\). The usual scale \(F^s_{\infty,q}(\mathbb R^n)\), with \(q<\infty\), is defined through local averages and will not be used in this chapter; see \cite[Section~1.5.2]{Triebel2}.

If $p=q<\infty$, Theorem~\ref{teo:tonelli} directly shows that
\[
B^s_{p,p}(\mathbb R^n)=F^s_{p,p}(\mathbb R^n)
\]
with the same norm for a fixed resolution. For distinct indices, interchanging $L^p(\mathbb R^n)$ and $\ell^q(\mathbb N_0)$ generally changes the space.

Reconstruction of a distribution requires retaining information about the frequencies of its summands. The next lemma expresses this condition and also allows treatment of the endpoint \(p=1\).

\begin{lemma}[Synthesis of blocks with controlled spectral support]
\label{lem:sintesis-bloques-espectrales}
Let \(s\in\mathbb R\), \(1\leq q\leq\infty\), and \(0<c<C\). Suppose that \(v_j\in L^p(\mathbb R^n)\), \(j\in\mathbb N_0\), and
\[
\operatorname{supp}\widehat v_0\subseteq\{\xi\in\mathbb R^n\mid\|\xi\|\leq C\},
\qquad
\operatorname{supp}\widehat v_j
\subseteq\{\xi\in\mathbb R^n\mid c2^j\leq\|\xi\|\leq C2^j\}
\quad(j\in\mathbb N).
\]
If \(1\leq p\leq\infty\), then
\[
\left\|\sum_{j\in\mathbb N_0}v_j\right\|_{B^s_{p,q}(\mathbb R^n)}
\leq C_1
\left\|\bigl(2^{js}\|v_j\|_{L^p(\mathbb R^n)}\bigr)_{j\in\mathbb N_0}\right\|_{\ell^q(\mathbb N_0)}.
\]
If \(1\leq p<\infty\), also
\[
\left\|\sum_{j\in\mathbb N_0}v_j\right\|_{F^s_{p,q}(\mathbb R^n)}
\leq C_2
\left\|\bigl(2^{js}v_j\bigr)_{j\in\mathbb N_0}\right\|_{L^p(\mathbb R^n;\ell^q(\mathbb N_0))}.
\]
In each assertion, finiteness of the right-hand side ensures convergence of the series in \(\mathcal S'(\mathbb R^n)\).
\end{lemma}

\begin{proof}
Choose multipliers \(P_j\) equal to one on the indicated supports, supported on slightly larger annuli when \(j\geq1\). For \(\zeta\in\mathcal S(\mathbb R^n)\) and any \(N>0\),
\[
\|P_j\zeta\|_{L^{p'}(\mathbb R^n)}
\leq C_N2^{-jN}\mathfrak p_N(\zeta),\qquad j\in\mathbb N_0,
\]
where \(\mathfrak p_N\) is a finite sum of Schwartz seminorms. To obtain this bound, differentiate the product of the symbol of \(P_j\) and \(\widehat\zeta\), use decay of all derivatives of \(\widehat\zeta\) on the annulus, and integrate by parts in the inversion formula. The number of integrations is chosen greater than \(n\); this simultaneously controls the uniform norm and the integral norm of the result, and the remaining norms follow from Hölder's inequality.

Each sequence norm in the statement dominates \(2^{js}\|v_j\|_{L^p(\mathbb R^n)}\). Since \(\langle v_j,\zeta\rangle=\langle v_j,P_j\zeta\rangle\), Hölder's inequality gives
\[
|\langle v_j,\zeta\rangle|
\leq C_N2^{-j(N+s)}\mathfrak p_N(\zeta)\,\mathcal N,
\]
where \(\mathcal N\) is the corresponding right-hand side. Choosing \(N>|s|\) proves absolute convergence of the pairing and its continuity on \(\zeta\).

The spectral supports imply \(\Delta_kv_j=0\) if \(|k-j|>N_0\), for a fixed integer \(N_0\). Young's inequality gives
\[
2^{ks}\left\|\Delta_k\sum_{j\in\mathbb N_0}v_j\right\|_{L^p(\mathbb R^n)}
\leq C\sum_{\substack{j\in\mathbb N_0\\|j-k|\leq N_0}}
2^{js}\|v_j\|_{L^p(\mathbb R^n)}.
\]
Proposition~\ref{prop:young-sucesiones} proves the Besov bound. For Triebel--Lizorkin, choose \(0<a<1\), so that \(a<p\) and, if \(q<\infty\), also \(a<q\). Lemma~\ref{lem:multiplicadores-diadicos-vectoriales}, applied to the blocks \(v_j\), allows the \(L^p(\mathbb R^n;\ell^q(\mathbb N_0))\) norm to be taken in the same finite sum of shifts. This gives the second bound. Apply the estimates first to finite sums; distributional convergence identifies each block of the limit, and Fatou's lemma allows passage to the series. For \(q=\infty\), use the supremum over \(j\in\mathbb N_0\).
\end{proof}

\begin{theorem}[Independence of the dyadic resolution]
\label{teo:independencia-resolucion-besov-triebel}
Let \(\|\cdot\|_{A,\Delta}\) and \(\|\cdot\|_{A,\nabla}\) be the norms constructed from two dyadic resolutions, where \(A=B^s_{p,q}(\mathbb R^n)\) or \(A=F^s_{p,q}(\mathbb R^n)\), with the indices in Definition~\ref{def:besov-triebel-lizorkin}. There are \(c_A,C_A>0\), independent of \(u\), such that
\[
c_A\|u\|_{A,\Delta}\leq\|u\|_{A,\nabla}\leq C_A\|u\|_{A,\Delta}.
\]
These spaces are Banach, and their inclusions into \(\mathcal S'(\mathbb R^n)\) are continuous.
\end{theorem}

\begin{proof}
Inserting the enlarged blocks gives
\[
\nabla_ku=
\sum_{\substack{j\in\mathbb N_0\\|j-k|\leq N_0}}
\nabla_k\widetilde\Delta_j\Delta_ju,\qquad k\in\mathbb N_0.
\]
The kernels have uniformly bounded integral norms. Young's inequalities for convolutions and sequences prove the Besov bound. For the Triebel--Lizorkin norm, Lemma~\ref{lem:multiplicadores-diadicos-vectoriales} applies to \(f_j=\Delta_ju\), whose spectral supports are controlled. The sum contains only \(2N_0+1\) index shifts and therefore gives the required bound. Interchanging the resolutions proves equivalence.

For completeness, take the sequence spaces
\[
b^s_{p,q}:=\left\{(v_j)_{j\in\mathbb N_0}\mid
\left\|\bigl(2^{js}\|v_j\|_{L^p(\mathbb R^n)}\bigr)_{j\in\mathbb N_0}\right\|_{\ell^q(\mathbb N_0)}
<\infty\right\},
\]
\[
f^s_{p,q}:=\left\{(v_j)_{j\in\mathbb N_0}\mid
\left\|\bigl(2^{js}v_j\bigr)_{j\in\mathbb N_0}\right\|_{L^p(\mathbb R^n;\ell^q(\mathbb N_0))}
<\infty\right\}.
\]
In \(L^p(\mathbb R^n;\ell^\infty(\mathbb N_0))\), we mean sequences of measurable functions, identified coordinatewise almost everywhere, with norm \(\displaystyle \|\displaystyle\sup_{j\in\mathbb N_0}|v_j|\|_{L^p(\mathbb R^n)}\). This is a mixed norm; the map taking values in \(\ell^\infty(\mathbb N_0)\) need not be strongly measurable in the Bochner sense. Completeness of both sequence spaces follows directly: from a Cauchy sequence, extract a subsequence whose differences have summable norms. The sum of the absolute values of these differences is finite almost everywhere in the corresponding sequence norm. This gives a coordinatewise measurable limit, and the summable-tail inequality proves convergence in norm. For the remaining indices, this agrees with the usual completeness of Lebesgue spaces and Banach sums. In each space, retain only the sequences whose Fourier supports lie in the annuli of the preceding lemma, chosen to contain the supports of \(\phi_j\). This subspace is closed: convergence in either norm implies convergence of each coordinate in \(L^p(\mathbb R^n)\) and hence in \(\mathcal S'\); a distribution vanishing outside a fixed closed set retains this property in the limit.

On this subspace, define
\[
R(v_j)_{j\in\mathbb N_0}:=\sum_{j\in\mathbb N_0}v_j,
\qquad
Su:=(\Delta_ju)_{j\in\mathbb N_0}.
\]
The lemma proves continuity of \(R\); the definitions of the norms prove continuity of \(S\). Dyadic reconstruction gives \(RS=I\). Thus \(SR\) is a continuous projection from the sequence subspace onto \(S(A)\). Its image is closed and Banach, and \(S\) identifies \(A\) isometrically with it. The pairing estimate obtained in the lemma, applied to \(v_j=\Delta_ju\), also proves continuity of \(A\hookrightarrow\mathcal S'\).
\end{proof}

For interpolation of Bessel spaces, it will be useful to have synthesis defined on sequences without spectral restrictions.

\begin{lemma}[Synthesis on the full sequence space]
\label{lem:sintesis-sucesiones-completas}
The operator
\[
\widetilde R(v_j)_{j\in\mathbb N_0}
:=\sum_{j\in\mathbb N_0}\widetilde\Delta_jv_j
\]
is continuous from \(b^s_{p,q}\) to \(B^s_{p,q}(\mathbb R^n)\) for \(1\leq p,q\leq\infty\). It is also continuous from \(f^s_{p,q}\) to \(F^s_{p,q}(\mathbb R^n)\) if \(1<p<\infty\) and \(1\leq q\leq\infty\). In both cases, \(\widetilde RS=I\).
\end{lemma}

\begin{proof}
The enlarged blocks have kernels dominated by
\[
Ch_j(x),\qquad
h_j(x):=2^{jn}(1+2^j\|x\|)^{-n-1}.
\]
The Besov bound follows from Young's inequality and the spectral synthesis lemma. For \(1<p<\infty\), \(1<q<\infty\), domination by the maximal operator and Fefferman--Stein control \(\|(\widetilde\Delta_jv_j)_{j\in\mathbb N_0}\|_{L^p(\mathbb R^n;\ell^q(\mathbb N_0))}\). If \(q=\infty\), use
\[
\sup_{j\in\mathbb N_0}|\widetilde\Delta_jv_j|
\leq C\mathcal M\left(\sup_{j\in\mathbb N_0}|v_j|\right).
\]
For \(q=1\), take a finite family and a function \(g\in L^{p'}(\mathbb R^n)\), \(g\geq0\). Tonelli's theorem and symmetry of \(h_j\) give
\[
\begin{aligned}
\int_{\mathbb R^n}g\sum_{j=0}^{N}|\widetilde\Delta_jv_j|
&\leq C\int_{\mathbb R^n}\sum_{j=0}^{N}|v_j|(h_j*g)\\
&\leq C\int_{\mathbb R^n}
\left(\sum_{j=0}^{N}|v_j|\right)\mathcal Mg\\
&\leq C_p
\left\|\sum_{j=0}^{N}|v_j|\right\|_{L^p(\mathbb R^n)}\|g\|_{L^{p'}(\mathbb R^n)}.
\end{aligned}
\]
Duality of \(L^p(\mathbb R^n)\) and monotone convergence yield the bound for the full family. Incorporate the weights \(2^{js}\) into \(v_j\). Then apply the spectral synthesis lemma to \(\widetilde\Delta_jv_j\).
\end{proof}

The requirement \(p>1\) in the second assertion matters: the proof of completeness for \(F^s_{1,q}(\mathbb R^n)\) uses the subspace with controlled spectral supports from the preceding theorem.

\begin{corollary}[Approximation and loss of regularity]
\label{cor:aproximacion-escalas-BF}
For \(p,q<\infty\), \(\mathcal S(\mathbb R^n)\) is dense in \(B^s_{p,q}(\mathbb R^n)\) and in \(F^s_{p,q}(\mathbb R^n)\). For any admissible fine index, \(1\leq p<\infty\), and \(\varepsilon>0\),
\[
A^s_{p,q}(\mathbb R^n)\hookrightarrow
B^{s-\varepsilon}_{p,p}(\mathbb R^n),\qquad A\in\{B,F\}.
\]
The dyadic partial sums converge in this latter space.
\end{corollary}

\begin{proof}
If both indices are finite, the tails of the block sequence tend to zero in the corresponding norm by dominated convergence. The spectral synthesis lemma gives the same conclusion for tails of the reconstruction. Each of the finitely many remaining blocks is approximated in \(L^p(\mathbb R^n)\) by Schwartz functions; applying a multiplier equal to one on its spectral support preserves convergence and keeps the support in a slightly larger annulus. Finite synthesis gives approximation in the norm of the space.

In every case,
\[
2^{js}\|\Delta_ju\|_{L^p(\mathbb R^n)}\leq\|u\|_{A^s_{p,q}},
\qquad j\in\mathbb N_0.
\]
Multiplying by \(2^{-j\varepsilon}\) and summing the \(p\)th powers proves the second assertion and convergence of the tails. For \(q=\infty\), the asserted convergence takes place in the space of order \(s-\varepsilon\).
\end{proof}

\subsection{Fundamental identifications and interpolation}

In this subsection, we use the abstract interpolation results developed in Chapter~\ref{cap:espacios-interpolacion}; in particular, we apply the real and complex methods to the Euclidean scales constructed here.

The newly defined scales contain several spaces introduced earlier. Two mechanisms should be distinguished. Identification of Bessel potentials with $F^s_{p,2}(\mathbb R^n)$ is a Littlewood--Paley statement and therefore depends on the quadratic combination of frequency blocks. In contrast, identification of noninteger Sobolev--Slobodeckij spaces with $B^s_{p,p}(\mathbb R^n)$ comes from real interpolation. The two proofs explain how different extensions of the integer Sobolev scale arise.

\begin{theorem}[Bessel potentials and integer-order Sobolev spaces]
\label{teo:identificaciones-H-F-W-entero}
Let $1<p<\infty$ and $s\in\mathbb R$.
\begin{enumerate}[label=(\alph*)]
\item For every $s\in\mathbb R$,
\[
H^{s,p}(\mathbb R^n)=F^s_{p,2}(\mathbb R^n)
\]
and there exist $c_{s,p},C_{s,p}>0$ such that
\[
c_{s,p}\|u\|_{H^{s,p}(\mathbb R^n)}
\leq\|u\|_{F^s_{p,2}(\mathbb R^n)}
\leq C_{s,p}\|u\|_{H^{s,p}(\mathbb R^n)}.
\]
\item If $m\in\mathbb N_0$, then
\[
H^{m,p}(\mathbb R^n)=W^{m,p}(\mathbb R^n)
\]
and there exist $c_{m,p},C_{m,p}>0$ such that
\[
c_{m,p}\|u\|_{W^{m,p}(\mathbb R^n)}
\leq\|u\|_{H^{m,p}(\mathbb R^n)}
\leq C_{m,p}\|u\|_{W^{m,p}(\mathbb R^n)}.
\]
\end{enumerate}
\end{theorem}

The classical identification $F^s_{p,2}(\mathbb R^n)=H^{s,p}(\mathbb R^n)$ can be found in \cite[Chapter~1]{Triebel2}.

\begin{proof}
First extend the Littlewood--Paley criterion to distributions. If \(v\in\mathcal S'(\mathbb R^n)\) and
\[
G_v:=\left(\sum_{j\in\mathbb N_0}|\Delta_jv|^2\right)^{\frac{1}{2}}
\in L^p(\mathbb R^n),
\]
dyadic reproduction and Cauchy--Schwarz give, for \(\zeta\in\mathcal S(\mathbb R^n)\),
\[
\begin{split}
|\langle v,\zeta\rangle|
&=\left|\sum_{j\in\mathbb N_0}
\langle\Delta_jv,\widetilde\Delta_j\zeta\rangle\right|\\
&\leq \|G_v\|_{L^p(\mathbb R^n)}
\left\|\left(\sum_{j\in\mathbb N_0}
|\widetilde\Delta_j\zeta|^2\right)^{\frac{1}{2}}\right\|_{L^{p'}(\mathbb R^n)}
\leq C\|G_v\|_{L^p(\mathbb R^n)}\|\zeta\|_{L^{p'}(\mathbb R^n)}.
\end{split}
\]
The sum of pairings converges absolutely by the same inequality. The argument with signs used in Theorem~\ref{teo:funcion-cuadrada-littlewood-paley} also applies to the enlarged blocks, since their symbols have the same uniform bounds and finite overlap. Density of \(\mathcal S\) in \(L^{p'}(\mathbb R^n)\) and duality of \(L^{p'}(\mathbb R^n)\) show that \(v\) is represented by a function in \(L^p(\mathbb R^n)\), with \(\|v\|_{L^p(\mathbb R^n)}\leq C\|G_v\|_{L^p(\mathbb R^n)}\). Together with the direct Littlewood--Paley bound, this proves the criterion in both directions for \(v\in\mathcal S'\).

On the support of \(\phi_j\), the weight \(\langle\xi\rangle^s\) is comparable to \(2^{js}\) uniformly in \(j\in\mathbb N_0\). To pass from this comparison of symbols to one of norms, consider
\[
m_j(\xi)=2^{-js}\langle\xi\rangle^s\widetilde\phi_j(\xi),
\qquad
\widetilde m_j(\xi)=2^{js}\langle\xi\rangle^{-s}
\widetilde\phi_j(\xi).
\]
After \(\xi=2^j\zeta\), all their derivatives are uniformly bounded on a fixed compact set. Integration by parts shows that their kernels satisfy the hypotheses of Lemma~\ref{lem:multiplicadores-diadicos-vectoriales}. Applying it with fine index \(2\), first to the blocks of \(u\) and then to those of \(J^su\), yields
\[
c\left\|\left(\sum_{j\in\mathbb N_0}
2^{2js}|\Delta_ju|^2\right)^{\frac{1}{2}}\right\|_{L^p(\mathbb R^n)}
\leq
\left\|\left(\sum_{j\in\mathbb N_0}
|\Delta_jJ^su|^2\right)^{\frac{1}{2}}\right\|_{L^p(\mathbb R^n)}
\leq
C\left\|\left(\sum_{j\in\mathbb N_0}
2^{2js}|\Delta_ju|^2\right)^{\frac{1}{2}}\right\|_{L^p(\mathbb R^n)}.
\]
The preceding distributional criterion states that these quantities are finite exactly when \(J^su\in L^p(\mathbb R^n)\), and are then equivalent to \(\|J^su\|_{L^p(\mathbb R^n)}\). This proves (a).

For (b), if \(\alpha\in\mathbb N_0^n\) and \(|\alpha|\leq m\), the symbol
\[
b_\alpha(\xi)=(i\xi)^\alpha\langle\xi\rangle^{-m}
\]
satisfies the Mikhlin conditions: each differentiation lowers its order by one, and its initial order is \(|\alpha|-m\leq0\). Thus \(D^\alpha u=T_{b_\alpha}J^mu\) and
\[
\sum_{\substack{\alpha\in\mathbb N_0^n\\|\alpha|\leq m}}
\|D^\alpha u\|_{L^p(\mathbb R^n)}
\leq C\|J^mu\|_{L^p(\mathbb R^n)}.
\]
To recover \(J^mu\) from the derivatives, introduce the positive polynomial
\[
P_m(\xi)=
\sum_{\substack{\alpha\in\mathbb N_0^n\\|\alpha|\leq m}}\xi^{2\alpha}.
\]
We have \(c_m\langle\xi\rangle^{2m}\leq P_m(\xi)
\leq C_m\langle\xi\rangle^{2m}\). The upper bound follows by estimating each monomial. For the lower bound, if \(\|\xi\|\geq1\), some coordinate satisfies \(|\xi_\ell|\geq\displaystyle\frac{\|\xi\|}{\sqrt n}\), and the summand \(\xi_\ell^{2m}\) gives the bound; if \(\|\xi\|\leq1\), use the constant term.

For the same multi-indices, define
\[
a_\alpha(\xi)=
\frac{\langle\xi\rangle^m(-i\xi)^\alpha}{P_m(\xi)}.
\]
Differentiating \(P_mP_m^{-1}=1\) and inducting on \(|\beta|\) gives
\[
|D^\beta(P_m^{-1})(\xi)|
\leq C_\beta\langle\xi\rangle^{-2m-|\beta|},
\qquad \beta\in\mathbb N_0^n.
\]
The Leibniz rule then implies \(|D^\beta a_\alpha(\xi)|
\leq C_{\alpha,\beta}\langle\xi\rangle^{|\alpha|-m-|\beta|}\); hence the \(a_\alpha\) are also Mikhlin symbols. The algebraic identity
\[
\sum_{\substack{\alpha\in\mathbb N_0^n\\|\alpha|\leq m}}
a_\alpha(\xi)(i\xi)^\alpha=\langle\xi\rangle^m
\]
gives, in \(\mathcal S'\),
\[
J^mu=
\sum_{\substack{\alpha\in\mathbb N_0^n\\|\alpha|\leq m}}
T_{a_\alpha}D^\alpha u.
\]
Continuity of each multiplier on \(L^p(\mathbb R^n)\) proves the reverse inequality. The sum norms and the \(\ell^p(\{\alpha\in\mathbb N_0^n\mid|\alpha|\leq m\})\) norms for this finite number of derivatives are equivalent, giving exactly the norm of \(W^{m,p}(\mathbb R^n)\) fixed in this book.
\end{proof}

The preceding theorem allows interpolation methods to be transferred to dyadic scales. For the real method, regularity is encoded in the weight $2^{js}$ of each block; for the complex method, these weights interpolate geometrically. This difference explains why the real method produces the Besov scale while the complex method preserves the Bessel scale.

\begin{theorem}[Interpolation of the scales]
\label{teo:interpolacion-escalas-BF}
Let $s_0<s_1$, $0<\theta<1$, and $s=(1-\theta)s_0+\theta s_1$.
\begin{enumerate}[label=(\alph*)]
\item If $1<p<\infty$ and $1\leq q\leq\infty$, then
\[
(H^{s_0,p}(\mathbb R^n),H^{s_1,p}(\mathbb R^n))_{\theta,q}
=
B^s_{p,q}(\mathbb R^n)
\]
and there exist $c,C>0$, depending on $n,s_0,s_1,p,q,\theta$ and the fixed resolution, such that
\[
c\|u\|_{B^s_{p,q}(\mathbb R^n)}
\leq\|u\|_{(H^{s_0,p}(\mathbb R^n),H^{s_1,p}(\mathbb R^n))_{\theta,q}}
\leq C\|u\|_{B^s_{p,q}(\mathbb R^n)}.
\]
\item If $1<p_0,p_1<\infty$ and
\[
\frac{1}{p}=\frac{1-\theta}{p_0}+\frac{\theta}{p_1},
\]
then
\[
[H^{s_0,p_0}(\mathbb R^n),H^{s_1,p_1}(\mathbb R^n)]_\theta
=
H^{s,p}(\mathbb R^n)
\]
and there exist $c,C>0$, depending on $n,s_0,s_1,p_0,p_1,\theta$ and the fixed resolution, such that
\[
c\|u\|_{H^{s,p}(\mathbb R^n)}
\leq\|u\|_{[H^{s_0,p_0}(\mathbb R^n),H^{s_1,p_1}(\mathbb R^n)]_\theta}
\leq C\|u\|_{H^{s,p}(\mathbb R^n)}.
\]
\item In particular, if $m\in\mathbb N$ and $0<\theta<1$, then
\[
(L^p(\mathbb R^n),W^{m,p}(\mathbb R^n))_{\theta,q}
=
B^{\theta m}_{p,q}(\mathbb R^n),
\]
and
\[
[L^p(\mathbb R^n),W^{m,p}(\mathbb R^n)]_\theta
=
H^{\theta m,p}(\mathbb R^n).
\]
\end{enumerate}
\end{theorem}

\begin{remark}[Why $p$ is fixed in the real method]
In part (a), $p$ cannot simply be replaced by two distinct exponents $p_0,p_1$. If $p^{-1}=(1-\theta)p_0^{-1}+\theta p_1^{-1}$, already in the low-frequency block one would obtain
\[
(L^{p_0}(\mathbb R^n),L^{p_1}(\mathbb R^n))_{\theta,q}
=
L^{p,q}(\mathbb R^n),
\]
a Lorentz space, rather than $L^p(\mathbb R^n)$ except in special cases. Formulas with $p_0\neq p_1$ lead to Besov--Lorentz scales. The complex method in part (b), in contrast, interpolates smoothness and the spatial exponent simultaneously within the Bessel scale.
\end{remark}

\begin{proof}
Define the analysis operator
\[
Su:=(\Delta_ju)_{j\in\mathbb N_0}
\]
and the synthesis operator
\[
R(v_j)_{j\in\mathbb N_0}
:=
\sum_{j\in\mathbb N_0}\widetilde\Delta_jv_j.
\]
Lemma~\ref{lem:sintesis-sucesiones-completas}, with $1<p<\infty$ and fine index $2$, proves that both operators are continuous on the scales below and satisfy $RS=I$. Thus the function scales are identified as retracts of the weighted sequence spaces associated with their dyadic blocks.

We prove (a). Write
\[
\delta:=s_1-s_0>0,
\qquad
t_k:=2^{-k\delta},
\qquad
a_j:=\|\Delta_ju\|_{L^p(\mathbb R^n)},
\quad j,k\in\mathbb N_0.
\]
For $k\in\mathbb N_0$, separate high and low frequencies by
\[
u=u_0^{(k)}+u_1^{(k)},
\]
with
\[
u_0^{(k)}:=\sum_{\substack{j\in\mathbb N_0\\j>k}}\widetilde\Delta_j\Delta_ju,
\qquad
u_1^{(k)}:=\sum_{j=0}^k\widetilde\Delta_j\Delta_ju.
\]
Finite interaction of the blocks, the triangle inequality, and the identification $H^{s_i,p}(\mathbb R^n)=F^{s_i}_{p,2}(\mathbb R^n)$ from Theorem~\ref{teo:identificaciones-H-F-W-entero} give
\begin{equation}
K(t_k,u;H^{s_0,p}(\mathbb R^n),H^{s_1,p}(\mathbb R^n))
\leq
C\left(
\sum_{\substack{j\in\mathbb N_0\\j>k}}2^{js_0}a_j
+t_k\sum_{j=0}^k2^{js_1}a_j
\right).
\label{eq:K-besov-cota-superior}
\end{equation}
In the reverse direction, let $u=v_0+v_1$ with $v_i\in H^{s_i,p}(\mathbb R^n)$. Uniform boundedness of the block $\Delta_j$ and the description $H^{s_i,p}(\mathbb R^n)=F^{s_i}_{p,2}(\mathbb R^n)$ imply
\[
2^{js_0}a_j
\leq
C\left(
\|v_0\|_{H^{s_0,p}(\mathbb R^n)}
+t_j\|v_1\|_{H^{s_1,p}(\mathbb R^n)}
\right).
\]
Taking the infimum over all decompositions of $u$ gives
\begin{equation}
2^{js}a_j
\leq
Ct_j^{-\theta}
K(t_j,u;H^{s_0,p}(\mathbb R^n),H^{s_1,p}(\mathbb R^n)).
\label{eq:K-besov-cota-inferior}
\end{equation}
Set $b_j:=2^{js}a_j$. Multiplying \eqref{eq:K-besov-cota-superior} by $t_k^{-\theta}$ transforms its right-hand side into
\[
C\left(
\sum_{\substack{j\in\mathbb N_0\\j>k}}2^{-\theta\delta(j-k)}b_j
+
\sum_{j=0}^k2^{-(1-\theta)\delta(k-j)}b_j
\right).
\]
Both coefficient sequences are geometrically summable. Lemma~\ref{lem:discretizacion-funcional-K} and the Hardy inequalities from Theorem~\ref{teo:desigualdades-hardy}, together with \eqref{eq:K-besov-cota-inferior}, give constants $c,C>0$ such that
\[
c\left\|
(2^{js}\|\Delta_ju\|_{L^p(\mathbb R^n)})_{j\in\mathbb N_0}
\right\|_{\ell^q(\mathbb N_0)}
\leq
\|u\|_{(H^{s_0,p}(\mathbb R^n),H^{s_1,p}(\mathbb R^n))_{\theta,q}}
\leq
C
\left\|
(2^{js}\|\Delta_ju\|_{L^p(\mathbb R^n)})_{j\in\mathbb N_0}
\right\|_{\ell^q(\mathbb N_0)}.
\]
This discretization describes the portion $0<t\leq1$. For $t\geq1$, the embedding $H^{s_1,p}(\mathbb R^n)\hookrightarrow H^{s_0,p}(\mathbb R^n)$ implies that there is $c_0>0$, independent of $t$ and $u$, such that
\[
c_0\|u\|_{H^{s_0,p}(\mathbb R^n)}
\leq
K(t,u;H^{s_0,p}(\mathbb R^n),H^{s_1,p}(\mathbb R^n))
\leq
\|u\|_{H^{s_0,p}(\mathbb R^n)}
\]
for every $t\geq1$. This norm is also controlled by the Besov norm. To see this, dyadic reconstruction and the triangle inequality give
\[
\|u\|_{H^{s_0,p}(\mathbb R^n)}
\leq
C\sum_{j\in\mathbb N_0}2^{js_0}a_j
=
C\sum_{j\in\mathbb N_0}2^{-j(s-s_0)}b_j.
\]
The sequence $(2^{-j(s-s_0)})_{j\in\mathbb N_0}$ belongs to $\ell^{q'}(\mathbb N_0)$ when $1<q<\infty$; for $q=1$ use its uniform bound, and for $q=\infty$ its summability. By Hölder's inequality for sequences, with these endpoint interpretations,
\[
\|u\|_{H^{s_0,p}(\mathbb R^n)}
\leq
C\|(b_j)_{j\in\mathbb N_0}\|_{\ell^q(\mathbb N_0)}.
\]
Thus the interpolation integral over $[1,\infty)$ is controlled by the same norm, since $t^{-\theta}\in L^q((1,\infty),\displaystyle\frac{dt}{t})$ for $q<\infty$ and is bounded for $q=\infty$. This proves (a), including the endpoints $q=1$ and $q=\infty$.

For (b), define, for $i\in\{0,1\}$,
\[
\ell^2_{s_i}(\mathbb N_0)
:=
\left\{(z_j)_{j\in\mathbb N_0}\;\middle|\;\sum_{j\in\mathbb N_0}2^{2js_i}|z_j|^2<\infty\right\},
\]
with its natural norm, and set
\[
Y_i:=L^{p_i}(\mathbb R^n;\ell^2_{s_i}(\mathbb N_0)).
\]
By Theorem~\ref{teo:identificaciones-H-F-W-entero}, $S$ and $R$ identify $H^{s_i,p_i}(\mathbb R^n)$ with a retract of $Y_i$. The weights interpolate geometrically,
\[
(2^{js_0})^{1-\theta}(2^{js_1})^\theta=2^{js},
\]
and the spatial exponent satisfies $p^{-1}=(1-\theta)p_0^{-1}+\theta p_1^{-1}$. Proposition~\ref{prop:interpolacion-ellp-valores-banach}, applied to weighted sequences on $\mathbb N_0$ (or equivalently, after extending them by zero to $\mathbb Z$), and then to the Bochner spaces, gives
\[
[Y_0,Y_1]_\theta
=
L^p(\mathbb R^n;\ell^2_s(\mathbb N_0)).
\]
Theorem~\ref{teo:interpolacion-compleja-operadores}, applied to $S$ and $R$, preserves the retraction. Finally, the identification
\[
H^{s,p}(\mathbb R^n)=F^s_{p,2}(\mathbb R^n)
\]
from Theorem~\ref{teo:identificaciones-H-F-W-entero} identifies the resulting retract with $H^{s,p}(\mathbb R^n)$ and proves (b).

For (c), use $H^{0,p}(\mathbb R^n)=L^p(\mathbb R^n)$ and $H^{m,p}(\mathbb R^n)=W^{m,p}(\mathbb R^n)$, both from Theorem~\ref{teo:identificaciones-H-F-W-entero}, and apply (a) and (b) with $s_0=0$ and $s_1=m$.
\end{proof}

\begin{corollary}[Real interpolation with different fine indices]
\label{cor:interpolacion-real-BF-indices-finos}
Let \(s_0<s_1\), \(0<\theta<1\), \(s=(1-\theta)s_0+\theta s_1\), \(1<p<\infty\), and \(1\leq q,q_0,q_1\leq\infty\). If \(A_0,A_1\in\{B,F\}\), then
\[
\bigl((A_0)^{s_0}_{p,q_0}(\mathbb R^n),
      (A_1)^{s_1}_{p,q_1}(\mathbb R^n)\bigr)_{\theta,q}
=B^s_{p,q}(\mathbb R^n)
\]
with equivalent norms. Here \((A_i)^{s_i}_{p,q_i}\) means \(B^{s_i}_{p,q_i}(\mathbb R^n)\) or \(F^{s_i}_{p,q_i}(\mathbb R^n)\) according to the choice of \(A_i\), for \(i\in\{0,1\}\).
\end{corollary}

\begin{proof}
The proof of part (a) of the preceding theorem uses two properties of the endpoints \(Y_i=H^{s_i,p}(\mathbb R^n)\):
\[
2^{js_i}\|\Delta_jv\|_{L^p(\mathbb R^n)}\leq C\|v\|_{Y_i},
\qquad
\left\|\sum_{j\in J}v_j\right\|_{Y_i}
\leq C\sum_{j\in J}2^{js_i}\|v_j\|_{L^p(\mathbb R^n)},
\]
where \(J\subseteq\mathbb N_0\) and the \(v_j\) have spectral support in the corresponding annuli. The second assertion includes series when the right-hand side is finite. Both properties also hold for \(Y_i=(A_i)^{s_i}_{p,q_i}\): the first follows from the definition, and the second from Lemma~\ref{lem:sintesis-bloques-espectrales} and the triangle inequality.

Thus, with \(t_k=2^{-k(s_1-s_0)}\) and \(a_j=\|\Delta_ju\|_{L^p(\mathbb R^n)}\), the same frequency splitting yields
\[
K(t_k,u;Y_0,Y_1)
\leq C\left(
\sum_{\substack{j\in\mathbb N_0\\j>k}}2^{js_0}a_j
+t_k\sum_{j=0}^k2^{js_1}a_j\right),
\qquad
2^{ks}a_k\leq Ct_k^{-\theta}K(t_k,u;Y_0,Y_1).
\]
Multiplying the first inequality by \(t_k^{-\theta}\) produces exactly the two geometric sequences from the preceding proof. Young's inequality and discretization of \(K\) compare the Besov norm with the portion of the interpolation norm corresponding to \(0<t\leq1\). If \(u\in B^s_{p,q}(\mathbb R^n)\), the second property, with \(i=0\), also gives
\[
\|u\|_{Y_0}\leq C
\sum_{j\in\mathbb N_0}2^{-j(s-s_0)}
\bigl(2^{js}a_j\bigr)
\leq C'\|u\|_{B^s_{p,q}(\mathbb R^n)}.
\]
For \(t\geq1\), use \(K(t,u;Y_0,Y_1)\leq\|u\|_{Y_0}\); integrability of \(t^{-\theta}\) completes the bound on that interval. This proves both inclusions.
\end{proof}

The real interpolation formula now identifies the spaces defined through differences with the Besov scale.

\begin{corollary}[Sobolev--Slobodeckij, Besov, and Triebel--Lizorkin]
\label{teo:identificaciones-H-W-B-F}
Let $1<p<\infty$ and $s>0$, and assume $s\notin\mathbb N$. Then
\[
W^{s,p}(\mathbb R^n)
=
B^s_{p,p}(\mathbb R^n)
=
F^s_{p,p}(\mathbb R^n)
\]
and there are constants $c_{n,s,p},C_{n,s,p}>0$ such that
\[
c_{n,s,p}\|u\|_{W^{s,p}(\mathbb R^n)}
\leq\|u\|_{B^s_{p,p}(\mathbb R^n)}=\|u\|_{F^s_{p,p}(\mathbb R^n)}
\leq C_{n,s,p}\|u\|_{W^{s,p}(\mathbb R^n)}.
\]
\end{corollary}

The identification through real interpolation is classical; see \cite[Chapters~1--2]{Triebel2}.

\begin{proof}
Write $s=m+\sigma$, with $m\in\mathbb N_0$ and $0<\sigma<1$. Corollary~\ref{cor:slobodeckij-interpolacion-altos-ordenes} gives
\[
W^{s,p}(\mathbb R^n)
=
(W^{m,p}(\mathbb R^n),W^{m+1,p}(\mathbb R^n))_{\sigma,p}.
\]
By Theorem~\ref{teo:identificaciones-H-F-W-entero}, the two integer-order spaces in the pair may be replaced by $H^{m,p}(\mathbb R^n)$ and $H^{m+1,p}(\mathbb R^n)$. Part (a) of Theorem~\ref{teo:interpolacion-escalas-BF} then gives
\[
W^{s,p}(\mathbb R^n)=B^s_{p,p}(\mathbb R^n).
\]
On the other hand, for a fixed dyadic resolution, Theorem~\ref{teo:tonelli} allows sum and integral to be interchanged:
\[
\begin{aligned}
\|u\|_{F^s_{p,p}(\mathbb R^n)}^p
&=
\int_{\mathbb R^n}\sum_{j\in\mathbb N_0}2^{jsp}|\Delta_ju(x)|^p\,dx\\
&=
\sum_{j\in\mathbb N_0}2^{jsp}\|\Delta_ju\|_{L^p(\mathbb R^n)}^p
=
\|u\|_{B^s_{p,p}(\mathbb R^n)}^p.
\end{aligned}
\]
This yields the identification and equivalence of norms.
\end{proof}

\begin{corollary}[When $H^{s,p}(\mathbb R^n)$ and $W^{s,p}(\mathbb R^n)$ coincide]
\label{cor:comparacion-H-W}
Let $s>0$ and $s\notin\mathbb N$, and let $1<p<\infty$.
\begin{enumerate}[label=(\alph*)]
\item If $1<p\leq2$, then $W^{s,p}(\mathbb R^n)\hookrightarrow H^{s,p}(\mathbb R^n)$.
\item If $2\leq p<\infty$, then $H^{s,p}(\mathbb R^n)\hookrightarrow W^{s,p}(\mathbb R^n)$.
\item Equality holds if $p=2$. If $p\neq2$, the corresponding embedding is strict.
\end{enumerate}
For $s=m\in\mathbb N_0$, however,
\[
H^{m,p}(\mathbb R^n)=W^{m,p}(\mathbb R^n)
\]
for every $1<p<\infty$.
\end{corollary}

\begin{proof}
By Corollary~\ref{teo:identificaciones-H-W-B-F} and Theorem~\ref{teo:identificaciones-H-F-W-entero}, it suffices to compare \(F^s_{p,p}(\mathbb R^n)\) with \(F^s_{p,2}(\mathbb R^n)\). The inclusion \(\ell^{q_0}(\mathbb N_0)\hookrightarrow
\ell^{q_1}(\mathbb N_0)\), for \(q_0\leq q_1\), applied at each point to \((2^{js}|\Delta_ju(x)|)_{j\in\mathbb N_0}\), proves the embeddings and equality for \(p=2\).

We construct examples establishing strictness inside any open ball \(B\subseteq\mathbb R^n\). This will also allow their use on domains. By independence of the resolution, we may choose a radial cutoff \(\chi\) equal to one when \(\|\xi\|\leq1\), zero when \(\|\xi\|\geq\displaystyle\frac{3}{2}\), and nonincreasing in the radius. For \(\phi_j(\xi)=\chi(2^{-j}\xi)-\chi(2^{-j+1}\xi)\), \(j\geq1\), the block equals one when \(\displaystyle\frac{3}{4}<2^{-j}\|\xi\|<1\). Take \(\xi_*=(\displaystyle\frac{7}{8},0,\ldots,0)\), a number \(0<\varepsilon<\displaystyle\frac{1}{16}\), and
\[
j_m=J+4m,\qquad
\lambda_m=2^{j_m}\xi_*,\qquad m\in\mathbb N,
\]
with \(J\in\mathbb N\). The balls of center \(\lambda_m\) and radius \(\varepsilon2^{j_m}\) lie in regions where \(\phi_{j_m}=1\); the other blocks vanish there.

Fix \(0\neq\zeta\in C_c^\infty(B)\) and a cutoff \(\omega\in C_c^\infty(B_{\mathrm{euc}}(0,\varepsilon))\) equal to one near the origin. Define
\[
\zeta_m=\mathcal F^{-1}
\bigl(\omega(2^{-j_m}\,\cdot)\widehat\zeta\bigr).
\]
For any multi-indices \(\alpha,\beta\in\mathbb N_0^n\) and any \(N>0\), rapid decay of \(\widehat\zeta\) gives
\begin{equation}
\sup_{x\in\mathbb R^n}
|x^\alpha D^\beta(\zeta_m-\zeta)(x)|
\leq C_{\alpha,\beta,N}2^{-j_mN}.
\label{eq:aproximacion-espectral-corte-lacunario}
\end{equation}
Indeed, the symbol \(1-\omega(2^{-j_m}\xi)\) vanishes on a ball of radius comparable to \(2^{j_m}\). Differentiating in \(\xi\) to produce \(x^\alpha\), multiplying by \(\xi^\beta\), and estimating the inverse transform by its \(L^1(\mathbb R^n)\) norm reduces the bound to tail integrals of a Schwartz function. These tails decay with every power of the radius.

Let \((a_m)_{m\in\mathbb N}\) be a bounded sequence. Since \(s>0\), the series
\[
w=\sum_{m\in\mathbb N}
a_m2^{-j_ms}e^{i\langle\lambda_m,x\rangle}\zeta_m(x),
\qquad
u=\zeta(x)\sum_{m\in\mathbb N}
a_m2^{-j_ms}e^{i\langle\lambda_m,x\rangle}
\]
converge uniformly and in \(\mathcal S'\). The second function has support in \(B\). By the choice of frequencies,
\[
2^{j_ms}\Delta_{j_m}w
=a_me^{i\langle\lambda_m,x\rangle}\zeta_m(x),
\qquad
\Delta_jw=0
\quad\text{if }j\notin\{j_m\mid m\in\mathbb N\}.
\]
The sequence \((a_m(\zeta_m-\zeta))_{m\in\mathbb N}\) belongs to \(L^p(\mathbb R^n;\ell^q(\mathbb N))\) for every \(1\leq q\leq\infty\): by \eqref{eq:aproximacion-espectral-corte-lacunario}, the sum of its \(L^p(\mathbb R^n)\) norms is finite. On the other hand,
\[
\bigl\|(a_m\zeta)_{m\in\mathbb N}\bigr\|_{L^p(\mathbb R^n;\ell^q(\mathbb N_0))}
=\|\zeta\|_{L^p(\mathbb R^n)}\|(a_m)_{m\in\mathbb N}\|_{\ell^q(\mathbb N)}.
\]
The triangle inequality, also applied to the difference, shows that
\[
w\in F^s_{p,q}(\mathbb R^n)
\quad\Longleftrightarrow\quad
(a_m)_{m\in\mathbb N}\in\ell^q(\mathbb N).
\]
Moreover, \(u-w\in\mathcal S\). To verify this, differentiate each summand of the difference. A derivative of the exponential produces at most a power of \(2^{j_m}\), and any polynomial weight in \(x\) is absorbed into the seminorm in \eqref{eq:aproximacion-espectral-corte-lacunario}. Choosing \(N\) larger than these powers makes the series summable in every Schwartz seminorm. Therefore,
\[
u\in F^s_{p,q}(\mathbb R^n)
\quad\Longleftrightarrow\quad
(a_m)_{m\in\mathbb N}\in\ell^q(\mathbb N).
\]

If \(p<2\), choose \(a_m=m^{-\gamma}\) with \(\displaystyle\frac{1}{2}<\gamma\leq\displaystyle\frac{1}{p}\); then \((a_m)\in\ell^2(\mathbb N_0)\setminus\ell^p(\mathbb N_0)\) and \(u\in H^{s,p}(\mathbb R^n)\setminus W^{s,p}(\mathbb R^n)\). If \(p>2\), choose \(\displaystyle\frac{1}{p}<\gamma\leq\displaystyle\frac{1}{2}\) to obtain \(u\in W^{s,p}(\mathbb R^n)\setminus H^{s,p}(\mathbb R^n)\). These choices prove strictness, even with compact support in \(B\). The integer-order case follows from Theorem~\ref{teo:identificaciones-H-F-W-entero}.
\end{proof}

\begin{remark}[Two different extensions of the integer scale]
The corollary clarifies a common source of terminological ambiguity. Some texts call both $H^{s,p}(\mathbb R^n)$ and $W^{s,p}(\mathbb R^n)$ ``fractional Sobolev spaces.'' We shall retain different letters. The scale
\[
H^{s,p}(\mathbb R^n)=F^s_{p,2}(\mathbb R^n)
\]
is stable under complex interpolation and suited to multipliers and elliptic operators. The scale
\[
W^{s,p}(\mathbb R^n)=B^s_{p,p}(\mathbb R^n)=F^s_{p,p}(\mathbb R^n),
\qquad s\notin\mathbb N,
\]
arises from real interpolation and is suited to differences, traces, and nonlocal energies. For $p=2$, both perspectives come together.
\end{remark}

\section{Multipliers, changes of coordinates, and localization}
\label{sec:localizacion-espacios-orden-real}

Euclidean theory can be transferred to a manifold only if two operations are stable: multiplication by a function in a partition of unity and change of coordinates. The next result treats both operations. Uniformity of the constants will be decisive in the chapter on bounded geometry.

\begin{definition}[The spaces \texorpdfstring{$BC^j$}{BCj} and \texorpdfstring{$BUC^j$}{BUCj}]
\label{def:espacios-bc-buc}\glsadd{espacios-bc-buc}
\index{BCj space@space $BC^j$}
\index{BUCj space@space $BUC^j$}
Let $U\subseteq\mathbb R^m$ be an open subset, let $F$ be a Banach space, and let $j\in\mathbb N_0$. Adopt the conventions $D^0f=f$ and $\mathcal L^0(\mathbb R^m;F)=F$. For $0\leq\ell\leq j$, the space $\mathcal L^\ell(\mathbb R^m;F)$ carries the norm
\[
 \|A\|_{\mathcal L^\ell}
 :=\sup\bigl\{\|A(v_1,\ldots,v_\ell)\|_F
 \mid |v_1|\leq1,\ldots,|v_\ell|\leq1\bigr\}.
\]
Define
\[
 BC^j(U,F)
 :=\left\{f\in C^j(U,F)\middle|
 \sum_{\ell=0}^j\sup_{x\in U}
 \|D^\ell f(x)\|_{\mathcal L^\ell}<+\infty\right\},
\]
with the norm
\[
 \|f\|_{BC^j(U,F)}
 :=\sum_{\ell=0}^j\sup_{x\in U}
 \|D^\ell f(x)\|_{\mathcal L^\ell}.
\]
The space $BUC^j(U,F)$ is the subspace of $BC^j(U,F)$ consisting of functions for which
\[
 D^\ell f\colon U\longrightarrow
 \mathcal L^\ell(\mathbb R^m;F)
\]
is uniformly continuous for each $0\leq\ell\leq j$; it carries the norm induced from $BC^j(U,F)$. In particular, $BC^0(U,F)$ is the space of continuous bounded functions, while $BUC^0(U,F)$ additionally requires uniform continuity.

Now let $Q_+^m=[0,1)\times(-1,1)^{m-1}$ and let $V\subseteq Q_+^m$ be a relatively open subset. Set
\[
 V^\circ:=V\cap\bigl((0,1)\times(-1,1)^{m-1}\bigr).
\]
By $C^j(V,F)$, we mean the space of functions $f\colon V\longrightarrow F$ for which $f\restriction_{V^\circ}$ is of class $C^j$ and each derivative $D^\ell(f\restriction_{V^\circ})$, $0\leq\ell\leq j$, admits a continuous extension to $V$, requiring the order-zero extension to agree with $f$. The extension is unique because $V^\circ$ is dense in $V$ and is also denoted by $D^\ell f$. Define $BC^j(V,F)$ and its norm by the preceding formulas with $V$ in place of $U$, and let $BUC^j(V,F)$ be the subspace for which each extended function
\[
 D^\ell f\colon V\longrightarrow
 \mathcal L^\ell(\mathbb R^m;F)
\]
is uniformly continuous with respect to the induced Euclidean distance. In particular, if $V=Q_+^m$, each derivative of a function in $BUC^j(Q_+^m,F)$ has a unique continuous extension to the Euclidean closure of $Q_+^m$. The convention also applies to overlap domains relatively open in a half-box and agrees with the usual notion of regularity up to the face $\{0\}\times(-1,1)^{m-1}$ on a manifold with boundary.

Finally, if $\Omega$ is any of the preceding domains,
\[
 BC^\infty(\Omega,F):=\bigcap_{j\in\mathbb N_0}BC^j(\Omega,F),
 \qquad
 BUC^\infty(\Omega,F):=\bigcap_{j\in\mathbb N_0}BUC^j(\Omega,F),
\]
with the families of seminorms $\|\cdot\|_{BC^j}$; no single uniform bound in $j$ is required. For matrix-valued functions and tensor coefficients, use the corresponding finite-dimensional space as codomain, with its operator norm or the tensor norm induced by the Euclidean norm. Thus the bounds hold simultaneously for all components.
\end{definition}

This is the sum convention used in Amann's uniform regularity. If one introduces the auxiliary norm
\(\displaystyle \|f\|_{\displaystyle\max,j}:=
 \displaystyle\max_{\ell\in\{0,\ldots,j\}}\displaystyle\sup_{x\in U}\|D^\ell f(x)\|_{\mathcal L^\ell},
\)
then explicitly,
\[
 \|f\|_{\max,j}\leq\|f\|_{BC^j}
 \leq(j+1)\|f\|_{\max,j}.
\]
The class of spaces is also independent of the chosen Euclidean norm. Indeed, if $a|v|_1\leq|v|_2\leq b|v|_1$ for every $v\in\mathbb R^m$, with $a,b>0$, then for each $\ell$,
\[
 \|A\|_{\mathcal L^\ell,2}
 \leq a^{-\ell}\|A\|_{\mathcal L^\ell,1},
 \qquad
 \|A\|_{\mathcal L^\ell,1}
 \leq b^\ell\|A\|_{\mathcal L^\ell,2}.
\]
Consequently,
\[
 \|f\|_{BC^j,2}
 \leq\max\{1,a^{-j}\}\|f\|_{BC^j,1},
 \qquad
 \|f\|_{BC^j,1}
 \leq\max\{1,b^j\}\|f\|_{BC^j,2}.
\]
The two domain metrics are bi-Lipschitz equivalent and therefore determine the same uniform continuity. If $F$ is finite dimensional and the codomain norms satisfy
\[
 c_F\|w\|_{F,1}\leq\|w\|_{F,2}\leq C_F\|w\|_{F,1}
 \qquad(w\in F),
\]
then the preceding inequalities become
\[
 \|f\|_{BC^j,2}
 \leq C_F\max\{1,a^{-j}\}\|f\|_{BC^j,1},
 \qquad
 \|f\|_{BC^j,1}
 \leq c_F^{-1}\max\{1,b^j\}\|f\|_{BC^j,2}.
\]
Thus the classes $BC^j$ and $BUC^j$ are independent of these choices.

Estimates for changes of coordinates require control of both frequency separation and spatial decay of the kernel. Write
\[
h_{r,D}(x):=2^{rn}(1+2^r\|x\|)^{-D},
\qquad r\in\mathbb N_0,\quad D>n,
\]
and, for a band-limited function \(v\),
\[
P_{r,\lambda}v(x):=
\sup_{z\in\mathbb R^n}
\frac{|v(z)|}{(1+2^r\|x-z\|)^\lambda}.
\]
Lemma~\ref{lem:control-maximal-peetre} ensures that \(P_{r,\lambda}v\leq C\mathcal M_a v\) when \(\lambda>\displaystyle\frac{n}{a}\) and the support of \(\widehat v\) lies in a ball of radius comparable to \(2^r\).

\begin{lemma}[Kernel estimate for a change of coordinates]
\label{lem:nucleo-cambio-coordenadas}
Let \(M\in\mathbb N\), let \(g\in BC^M(\mathbb R^n)\), and let \(\psi\) be a global diffeomorphism of class \(C^{M+1}\). Suppose that \(\psi\) and \(\psi^{-1}\) are Lipschitz and that their derivatives of orders \(1,\ldots,M+1\) are bounded. The kernel \(L_{j,r}(x,z)\) of
\[
\Delta_j\bigl[g\,((\widetilde\Delta_r v)\circ\psi)\bigr](x)
=\int_{\mathbb R^n}L_{j,r}(x,z)v(z)\,dz
\]
satisfies, for any \(D>n\),
\begin{equation}
|L_{j,r}(x,z)|
\leq C_{M,D}2^{-M|j-r|}
h_{\min\{j,r\},D}(\psi(x)-z),
\qquad j,r\in\mathbb N_0.
\label{eq:nucleo-cambio-diadico}
\end{equation}
\end{lemma}

\begin{proof}
Denote by \(K_j\) and \(\widetilde K_r\) the real convolution kernels of the blocks, including their Fourier constants. Then
\[
L_{j,r}(x,z)=
\int_{\mathbb R^n}K_j(x-y)g(y)\widetilde K_r(\psi(y)-z)\,dy.
\]
Positive-index kernels have all moments zero. Moreover, for any multi-index \(\alpha\) and any \(D_1>0\),
\[
|D^\alpha K_j(y)|+|D^\alpha\widetilde K_j(y)|
\leq C_{\alpha,D_1}
2^{j(n+|\alpha|)}(1+2^j\|y\|)^{-D_1}.
\]

If \(j\geq r\) and \(j\geq1\), subtract from \(y\mapsto g(y)\widetilde K_r(\psi(y)-z)\) its Taylor polynomial at \(x\) of degree \(M-1\). The polynomial terms vanish upon integration against \(K_j(x-y)\). The Leibniz and Faà di Bruno formulas, Theorem~\ref{faa di bruno multivariable}, show that derivatives of order \(M\) of the remainder are bounded by
\[
C2^{r(n+M)}
(1+2^r\|\psi(x)-z\|)^{-D}
(1+2^r\|x-y\|)^D
\]
on the segment between \(x\) and \(y\). Here we use
\[
1+2^r\|\psi(x)-z\|
\leq C(1+2^r\|\psi(x+t(y-x))-z\|)
(1+2^r\|x-y\|),\qquad 0\leq t\leq1.
\]
The integral remainder also contains \(\|x-y\|^M\). The substitution \(w=2^j(x-y)\), with \(D_1>M+D+n\), gives
\[
\int_{\mathbb R^n}|K_j(x-y)|\,\|x-y\|^M
(1+2^r\|x-y\|)^D\,dy\leq C2^{-jM}.
\]
This yields \eqref{eq:nucleo-cambio-diadico} in this case.

If \(r>j\), first make the substitution \(w=\psi(y)\). The factor against which \(\widetilde K_r(w-z)\) acts is
\[
K_j(x-\psi^{-1}(w))\,g(\psi^{-1}(w))\,|\det D\psi^{-1}(w)|.
\]
Subtract its Taylor polynomial at \(z\) of degree \(M-1\). Its derivatives of order \(M\) require at most \(M+1\) derivatives of \(\psi^{-1}\) because of the determinant. The same computation, now with \(2^{jM}\|w-z\|^M\), gives \(2^{-M(r-j)}\). The bi-Lipschitz condition compares \(\|x-\psi^{-1}(z)\|\) with \(\|\psi(x)-z\|\). When \(j=r=0\), the bound follows directly from the integral of two rapidly decreasing kernels. This also covers equal positive indices.

\end{proof}

\begin{lemma}[Multipliers and diffeomorphisms on the Sobolev scale]
\label{lema: multiplicadores y difeomorfismos sobolev}
Let \(s\in\mathbb R\) and \(k\in\mathbb N\), with \(k>|s|+n+1\). For \(H^{s,p}(\mathbb R^n)\), assume \(1<p<\infty\); for \(B^s_{p,q}(\mathbb R^n)\), assume \(1\leq p,q\leq\infty\); and for \(F^s_{p,q}(\mathbb R^n)\), assume \(1\leq p<\infty\), \(1\leq q\leq\infty\).
\begin{enumerate}[label=(\alph*)]
\item If \(g\in BC^k(\mathbb R^n)\), then
\[
\|gu\|_{A^s_{p,q}}\leq C\|g\|_{BC^k}\|u\|_{A^s_{p,q}},
\qquad A\in\{B,F\},
\]
and the same assertion holds in \(H^{s,p}(\mathbb R^n)\).
\item Let \(\psi\colon\mathbb R^n\longrightarrow\mathbb R^n\) be a diffeomorphism of class \(C^{k+1}\), with
\[
|\det D\psi(x)|\geq c>0,\qquad
\|D^\alpha\psi\|_{L^\infty(\mathbb R^n,\mathbb R^n)}\leq C_\alpha
\quad(1\leq|\alpha|\leq k+1).
\]
Then \(u\mapsto u\circ\psi\) is continuous on the same scales. The constants depend only on the parameters and the indicated bounds.
\item Let \(\psi\colon U\longrightarrow V\) be a diffeomorphism and let \(\eta\in C_c^{k+1}(U)\). Suppose that the derivative bounds and bi-Lipschitz conditions of the preceding lemma hold on a neighborhood of the support of \(\eta\). The operator
\[
T_{\eta,\psi}u(x):=
\begin{cases}
\eta(x)u(\psi(x)),&x\in U,\\
0,&x\notin U
\end{cases}
\]
is continuous on the same scales. The bound is uniform for families with the same derivative bounds, bi-Lipschitz constants, an upper bound for the diameters of the supports, and a positive lower bound for their distance from the boundary of these neighborhoods.
\end{enumerate}
\end{lemma}

\begin{proof}
The inverse matrix formula and the lower bound for the determinant control \(D\psi^{-1}\). Differentiating \(\psi^{-1}\circ\psi=I\) and applying Faà di Bruno, the term containing \(D^m\psi^{-1}\) is
\[
(D^m\psi^{-1})\circ\psi\,[D\psi,\ldots,D\psi].
\]
The other terms contain derivatives of \(\psi^{-1}\) of order less than \(m\). Induction on \(m=2,\ldots,k+1\), followed by multiplication by \(D\psi^{-1}\) in each argument, controls all derivatives of the inverse. Integrating the first derivative along segments proves that both functions are globally Lipschitz.

Choose \(0<a<1\), \(\lambda>\displaystyle\frac{n}{a}\), and an integer \(M<k\) so that \(M-\lambda>|s|\). This is possible because \(k>|s|+n+1\). Apply \eqref{eq:nucleo-cambio-diadico} to \(v_r=\Delta_ru\), with \(D>n+\lambda\). Peetre's definition gives
\[
|v_r(z)|
\leq P_{r,\lambda}v_r(\psi(x))
(1+2^r\|\psi(x)-z\|)^\lambda.
\]
The integral of the second factor against \(\displaystyle h_{\min\{j,r\},D}\) is bounded by \(C2^{\lambda(r-j)_+}\). Thus,
\begin{equation}
|\Delta_j[g\,((\widetilde\Delta_r\Delta_ru)\circ\psi)](x)|
\leq C2^{-(M-\lambda)|j-r|}
\mathcal M_a(\Delta_ru)(\psi(x)).
\label{eq:casi-diagonal-multiplicacion-composicion}
\end{equation}
Multiplication corresponds to \(\psi=I\), and composition to \(g=1\).

Multiply by \(2^{js}\) and sum in \(r\in\mathbb N_0\). The coefficients
\[
\bigl(2^{-(M-\lambda-|s|)|j-r|}\bigr)_{j-r\in\mathbb Z}
\]
form a summable sequence. For Besov spaces, scalar boundedness of the maximal operator on \(L^{\displaystyle\frac{p}{a}}(\mathbb R^n)\), followed by Young's inequality for sequences, gives the bound; if \(p=\infty\), use its bound on \(L^\infty(\mathbb R^n)\) directly. For Triebel--Lizorkin spaces, apply Fefferman--Stein with exponents \(\displaystyle\frac{p}{a}>1\) and \(\displaystyle\frac{q}{a}>1\). If \(q=\infty\), use the maximal operator of the supremum, as in Lemma~\ref{lem:multiplicadores-diadicos-vectoriales}. Change of variables introduces only the bound for \(|\det D\psi^{-1}|\). The identity \(H^{s,p}(\mathbb R^n)=F^s_{p,2}(\mathbb R^n)\) gives the assertion for Bessel spaces.

For part (c), reduce the local change to finitely many global changes. Fix \(x_0\) near the support of \(\eta\) and write
\[
A_0=D\psi(x_0),\qquad
a_0(x)=\psi(x_0)+A_0(x-x_0).
\]
If \(B(x_0,2r)\) lies in the neighborhood under consideration, choose a cutoff \(\chi\) equal to one on \(B(x_0,r)\) and supported in \(B(x_0,2r)\), with \(\|D\chi\|_{L^\infty(\mathbb R^n)}\leq \displaystyle\frac{C}{r}\). The fundamental theorem of calculus, applied to \(\psi-a_0\) and its derivative, gives
\[
\|\psi-a_0\|_{L^\infty(B(x_0,2r))}
\leq Cr^2\|D^2\psi\|_{L^\infty(B_{\mathrm{euc}}(x_0,2r))},\qquad
\|D(\psi-a_0)\|_{L^\infty(B(x_0,2r))}
\leq Cr\|D^2\psi\|_{L^\infty(B_{\mathrm{euc}}(x_0,2r))}.
\]
Extend \(F_0=\chi(\psi-a_0)\) by zero. By the Leibniz rule, \(\|DF_0\|_{L^\infty(B_{\mathrm{euc}}(x_0,2r))}\leq Cr\|D^2\psi\|_{L^\infty(B_{\mathrm{euc}}(x_0,2r))}\). Choose \(r>0\), uniformly under the assumptions of this part, so that \(\|A_0^{-1}\|\|DF_0\|_{L^\infty(\mathbb R^n)}<\displaystyle\frac{1}{2}\). The function
\[
\Psi_0(x)=a_0(x)+F_0(x),\qquad x\in\mathbb R^n,
\]
agrees with \(\psi\) on \(B(x_0,r)\) and is bi-Lipschitz. To check surjectivity as well, given \(y\in\mathbb R^n\), apply the fixed-point theorem to
\[
x\longmapsto x_0+A_0^{-1}
\bigl(y-\psi(x_0)-F_0(x)\bigr),
\]
whose contraction constant is less than \(\displaystyle\frac{1}{2}\). Its fixed point is the unique solution of \(\Psi_0(x)=y\). The inverse function theorem and the derivative bounds already established give a global diffeomorphism with the bounds required in (b).

Cover the support of \(\eta\) by finitely many of these balls and choose a subordinate smooth partition \((\beta_\ell)_{\ell=1}^N\). Then
\[
T_{\eta,\psi}u
=\sum_{\ell=1}^N\eta\beta_\ell\,(u\circ\Psi_\ell).
\]
Applying (a) and (b) to each summand proves (c). The ball radius, the number \(N\), and the cutoff derivatives are controlled by the bounds in the statement, giving the asserted uniformity as well.

Let us specify the distributional meaning of operators with coefficients of finite regularity for every $s\in\mathbb R$. Choose an integer $m$ such that $|s|+n<m<k$ and set
\[
Q_m(\zeta):=\sum_{|\alpha|\leq m}
\bigl(\|D^\alpha\zeta\|_{L^1(\mathbb R^n)}
     +\|D^\alpha\zeta\|_{L^\infty(\mathbb R^n)}\bigr).
\]
In the proof of Lemma~\ref{lem:sintesis-bloques-espectrales}, the reproducing multipliers $P_j$ can be chosen with even symbols. For $j\geq1$, their kernels have vanishing moments. Subtracting the Taylor polynomial of $\zeta(x-y)$ of degree $m-1$, using the integral remainder, and applying Minkowski gives
\[
\|P_j\zeta\|_{L^1}+\|P_j\zeta\|_{L^\infty}
\leq C_m2^{-jm}\sum_{|\alpha|=m}
\bigl(\|D^\alpha\zeta\|_{L^1}+\|D^\alpha\zeta\|_{L^\infty}\bigr).
\]
Indeed, each remainder norm is bounded by $C_m\|y\|^m$ times the sum of the norms of the derivatives of order $m$, and the absolute moment of the kernel at scale $2^{-j}$ is $C_m2^{-jm}$. For $j=0$, use Young. Hölder between $L^1$ and $L^\infty$ gives the same bound in $L^{p'}$. Since
$\|\Delta_ju\|_{L^p}\leq2^{-js}\|u\|_{A^s_{p,q}}$ and $m+s>0$, we obtain
\[
|\langle u,\zeta\rangle|
\leq C\|u\|_{A^s_{p,q}}Q_m(\zeta)
\sum_{j=0}^\infty2^{-j(m+s)}
=\frac{C}{1-2^{-(m+s)}}\|u\|_{A^s_{p,q}}Q_m(\zeta).
\]
Thus $u$ acts continuously on the closure of $\mathcal S$ for $Q_m$. This closure contains every $C_c^m$ function and every $C^m$ function whose derivatives of order at most $m$ are integrable and vanish at infinity. To verify this, first multiply by $\eta_R$: the Leibniz rule, integral tails, and uniform vanishing at infinity give convergence in $Q_m$. Then mollify the cutoff; translations are continuous in both $L^1$ and $C_0$, and mollification converges in each norm defining $Q_m$.

Define the operators by
\[
\langle gu,\zeta\rangle=\langle u,g\zeta\rangle,\qquad
\langle u\circ\psi,\zeta\rangle
=\left\langle u,(\zeta\circ\psi^{-1})|\det D\psi^{-1}|\right\rangle.
\]
For $\zeta\in\mathcal S$, both functions on the right belong to this closure. The Leibniz rule controls the product using the bounded derivatives of $g$; the chain rule, the bounds for $\psi^{-1}$, and change of variables control the second test function. The bi-Lipschitz condition also implies that $\psi^{-1}(x)$ escapes to infinity as $x$ does. The required derivatives have orders at most $m$ for $g$ and $m+1$ for $\psi^{-1}$, all of which are available. These pairings therefore define tempered distributions and agree with the usual operations on smooth functions.

Set $u_N=\sum_{j=0}^N\Delta_ju$. The synthesis lemma gives $\sup_N\|u_N\|_{A^s_{p,q}}\leq C_0\|u\|_{A^s_{p,q}}$, and Littlewood--Paley gives $u_N\to u$ in $\mathcal S'$. If $w$ is one of the preceding test functions of finite regularity and $w_\ell\in\mathcal S$ converges to $w$ in $Q_m$, the pairing bound gives
\[
|\langle u_N-u,w\rangle|
\leq C_1\|u\|_{A^s_{p,q}}Q_m(w-w_\ell)
     +|\langle u_N-u,w_\ell\rangle|.
\]
First choose $\ell$ to make the first term arbitrarily small; then $N\to\infty$ makes the second term tend to zero. Hence the operators on $u_N$ converge to the distributional operator defined above. The almost diagonal bounds already obtained apply to each finite sum. Tempered convergence identifies each block of the limit pointwise, and Fatou gives the final bounds. This argument includes $s=0$ and $q=\infty$ without requiring density in the norm of order $s$.
\end{proof}

To formulate the localization principle, call a family $(\eta_\nu)_{\nu\in I}\subseteq C_c^\infty(\mathbb R^n)$ a \textit{uniform localization system} if it satisfies:
\begin{enumerate}[label=(\roman*)]
\item $0\leq\eta_\nu\leq1$ and $\displaystyle\sum_{\nu\in I}\eta_\nu=1$;
\item each point belongs to at most $N$ supports;
\item support diameters and all $\|\eta_\nu\|_{BC^k(\mathbb R^n)}$ norms are uniformly bounded;
\item there exist $\chi_\nu\in C_c^\infty(\mathbb R^n)$, with $\chi_\nu=1$ near $\operatorname{supp}(\eta_\nu)$, satisfying the same bounds and the same type of overlap.
\end{enumerate}

\begin{theorem}[Localization principle]
\label{teo:principio-localizacion-F-H-W}
Let $(\eta_\nu)$ be a uniform localization system.
\begin{enumerate}[label=(\alph*)]
\item If $s\in\mathbb R$, $1<p<\infty$, and $1\leq q\leq\infty$, there are constants $c,C>0$, independent of $u$, such that
\[
c\|u\|_{F^s_{p,q}(\mathbb R^n)}
\leq
\left(
\displaystyle\sum_{\nu\in I}
\|\eta_\nu u\|_{F^s_{p,q}(\mathbb R^n)}^p
\right)^{\frac{1}{p}}
\leq
C\|u\|_{F^s_{p,q}(\mathbb R^n)}.
\]
In particular, the formula holds for $H^{s,p}(\mathbb R^n)=F^s_{p,2}(\mathbb R^n)$ and, if $s>0$ is not an integer, for $W^{s,p}(\mathbb R^n)=F^s_{p,p}(\mathbb R^n)$.
\item For $B^s_{p,q}(\mathbb R^n)$ with $q\neq p$, there is in general no analogous equivalence with an outer $\ell^p(I)$ sum. Such an equivalence does hold when $q=p$, because $B^s_{p,p}(\mathbb R^n)=F^s_{p,p}(\mathbb R^n)$.
\end{enumerate}
\end{theorem}

\begin{proof}
Omit indices for which \(\eta_\nu=0\) and choose \(x_\nu\in\operatorname{supp}\eta_\nu\), \(\nu\in I\). The diameter bounds place the supports of \(\eta_\nu\) and \(\chi_\nu\) in balls of fixed radius around \(x_\nu\). Moreover, there is \(r_0>0\), independent of \(\nu\), such that \(\chi_\nu>\displaystyle\frac{1}{2}\) on \(B(x_\nu,r_0)\): simply take
\[
r_0<\frac{1}{2\left(1+\sup_{\nu\in I}\|D\chi_\nu\|_{L^\infty(\mathbb R^n)}\right)}.
\]
Since the supports of \(\chi_\nu\) have bounded overlap, comparing the volumes of these balls gives
\[
\#\{\nu\in I\mid\|x-x_\nu\|\leq R\}\leq C(1+R)^n,
\qquad x\in\mathbb R^n,\quad R>0.
\]
In particular, \(I\) is countable and, for \(D>n\),
\begin{equation}
\sup_{x\in\mathbb R^n}
\sum_{\nu\in I}(1+\|x-x_\nu\|)^{-D}<\infty.
\label{eq:suma-pesos-localizacion}
\end{equation}
To verify the last bound, group the centers in \(\|x-x_\nu\|\leq1\) and in the annuli \(2^m<\|x-x_\nu\|\leq2^{m+1}\), \(m\in\mathbb N_0\). The resulting series is dominated by \(\displaystyle \sum_{m\in\mathbb N_0}2^{m(n-D)}\).

We shall use the preceding kernel estimate with the spatial decay supplied by the cutoffs. For \(a_\nu=\eta_\nu\) or \(a_\nu=\chi_\nu\), let \(L_{\nu,j,r}\) be the kernels of \(\Delta_j a_\nu\widetilde\Delta_r\). For any \(M,D,E>0\), with \(E>n\),
\begin{equation}
|L_{\nu,j,r}(x,z)|
\leq C_{M,D,E}2^{-M|j-r|}
(1+\|x-x_\nu\|)^{-D}
h_{\min\{j,r\},E}(x-z).
\label{eq:nucleo-localizacion-espacio-frecuencia}
\end{equation}
The constant is uniform in \(\nu,j,r\).

We explain the spatial factor. In the integral
\[
L_{\nu,j,r}(x,z)=
\int_{\mathbb R^n}K_j(x-y)a_\nu(y)\widetilde K_r(y-z)\,dy,
\]
\(y\) belongs to a ball of fixed radius around \(x_\nu\). If \(x\) lies in a slightly larger ball, that factor is comparable to one, and Lemma~\ref{lem:nucleo-cambio-coordenadas} applies. On the complement, \(1+\|x-x_\nu\|\leq C(1+\|x-y\|)\). For \(j\geq r\), decay of \(K_j(x-y)\) allows extraction of \((1+\|x-x_\nu\|)^{-D}2^{-M(j-r)}\); the remaining powers are integrated with \(\widetilde K_r\). For \(r>j\), subtract Taylor polynomials against \(\widetilde K_r\). Each derivative of \(K_j(x-y)a_\nu(y)\) retains the same decay outside the ball and contributes at most a factor \(2^j\). Integration of the remainder contributes \(2^{-rM}\). In both cases, the inequality
\[
1+2^{\min\{j,r\}}\|x-z\|
\leq
(1+2^{\min\{j,r\}}\|x-y\|)
(1+2^{\min\{j,r\}}\|y-z\|)
\]
allows the factor \(\displaystyle h_{\min\{j,r\},E}(x-z)\) to be retained. Initially choose decay powers larger than the sum of the extracted powers plus \(n+1\), so all remaining integrals are uniformly finite. Only finitely many derivative bounds of the cutoffs are involved.

Fix \(0<a<1\), \(\lambda>\displaystyle\frac{n}{a}\), \(D>n\), and \(M>|s|+\lambda\). Integrating \eqref{eq:nucleo-localizacion-espacio-frecuencia} against \(\Delta_ru\), as in \eqref{eq:casi-diagonal-multiplicacion-composicion}, gives
\[
2^{js}|\Delta_j(a_\nu u)(x)|
\leq Cw_\nu(x)
\sum_{r\in\mathbb N_0}2^{-(M-\lambda-|s|)|j-r|}
2^{rs}P_{r,\lambda}(\Delta_ru)(x),
\]
where \(w_\nu(x):=(1+\|x-x_\nu\|)^{-D}\). Take the \(\ell^q(\mathbb N_0)\) norm in \(j\), use Young's inequality for sequences, and set
\[
G_u(x):=
\left\|\bigl(2^{rs}P_{r,\lambda}(\Delta_ru)(x)\bigr)_{r\in\mathbb N_0}
\right\|_{\ell^q(\mathbb N_0)}.
\]
Peetre control and Fefferman--Stein, with the supremum argument if \(q=\infty\), imply
\[
\|G_u\|_{L^p(\mathbb R^n)}\leq C\|u\|_{F^s_{p,q}(\mathbb R^n)}.
\]
Thus \eqref{eq:suma-pesos-localizacion} and Tonelli's theorem give
\[
\sum_{\nu\in I}\|\eta_\nu u\|_{F^s_{p,q}(\mathbb R^n)}^p
\leq C\int_{\mathbb R^n}\sum_{\nu\in I}w_\nu(x)^pG_u(x)^p\,dx
\leq C\|u\|_{F^s_{p,q}(\mathbb R^n)}^p.
\]

For the reverse inequality, first consider a finite family \((v_\nu)_{\nu\in I}\) and the operator
\[
R(v_\nu)_{\nu\in I}:=\sum_{\nu\in I}\chi_\nu v_\nu.
\]
The same estimate, followed by the triangle inequality in \(\ell^q(\mathbb N_0)\), gives pointwise
\[
\left\|\bigl(2^{js}\Delta_jR(v_\nu)(x)\bigr)_{j\in\mathbb N_0}
\right\|_{\ell^q(\mathbb N_0)}
\leq C\sum_{\nu\in I}w_\nu(x)G_{v_\nu}(x).
\]
Hölder's inequality for the weighted sum over \(I\) implies
\[
\left(\sum_{\nu\in I}w_\nu G_{v_\nu}\right)^p
\leq
\left(\sum_{\nu\in I}w_\nu\right)^{p-1}
\sum_{\nu\in I}w_\nu G_{v_\nu}^p.
\]
Integration yields
\[
\|R(v_\nu)\|_{F^s_{p,q}(\mathbb R^n)}^p
\leq C\sum_{\nu\in I}\|v_\nu\|_{F^s_{p,q}(\mathbb R^n)}^p.
\]
Finite families are dense in the outer \(\ell^p(I)\) sum, so \(R\) extends continuously to it. Apply this bound to \(v_\nu=\eta_\nu u\). The sum is locally finite in the distributional sense, and \(\chi_\nu\eta_\nu=\eta_\nu\), so \(R(\eta_\nu u)=u\). This completes (a).

For (b), construct functions with compact spatial support. Choose a smooth partition obtained by lattice translation that equals one on small balls around its centers, and \(\zeta\in C_c^\infty(\mathbb R^n)\), \(\zeta\neq0\), supported in one of these balls. Choose distinct lattice centers \(x_m\) sufficiently far apart that each translate of \(\zeta\) lies in a region where only one cutoff equals one.

We may use a dyadic resolution whose annular symbol equals one near some point \(\xi_*\neq0\). Indeed, choose the initial cutoff equal to one on \(\|\xi\|\leq1\) and zero for \(\|\xi\|\geq\displaystyle\frac{3}{2}\); its dyadic difference equals one on \(\displaystyle\frac{3}{4}<\|\xi\|<1\). Fix \(\xi_*\) in that annulus. For integers \(j_m=J+Lm\), \(m\in\mathbb N\), set
\[
v_m(x):=
2^{-j_ms}\zeta(x-x_m)e^{i2^{j_m}\xi_*\cdot x},
\qquad
u_N:=\sum_{m=1}^{N}v_m.
\]
The integers \(J,L\) will be fixed once, independently of \(N\). For any \(M>|s|\), regularity of \(\zeta\) gives
\begin{equation}
2^{ks}\|\Delta_kv_m\|_{L^p(\mathbb R^n)}
\leq C_M2^{-(M-|s|)|k-j_m|},
\qquad k\in\mathbb N_0,\quad m\in\mathbb N.
\label{eq:bloques-modulacion-compacta}
\end{equation}
For \(k>j_m\), use \(M\) derivatives of \(\zeta e^{i2^{j_m}\xi_*\cdot x}\) and the moments of \(K_k\). For \(k<j_m\) separated by more than a fixed number of indices, integrate by parts in the direction \(\xi_*\) in the convolution formula: each integration contributes \(2^{-j_m}\), while differentiating the kernel contributes at most \(2^k\). Derivatives of \(\zeta\) have fixed \(L^p(\mathbb R^n)\) norms. The remaining indices are controlled by Young's inequality. These operations prove \eqref{eq:bloques-modulacion-compacta} also for \(p=1\) and \(p=\infty\).

Moreover,
\[
2^{j_ms}e^{-i2^{j_m}\xi_*\cdot x}\Delta_{j_m}v_m(x)
\longrightarrow\zeta(x-x_m)
\quad\text{in }L^p(\mathbb R^n)
\]
uniformly in the translations as \(j_m\to\infty\). To see this, the modulated kernel has mass one and scale \(2^{-j_m}\); write the difference using \(\zeta(\,\cdot-y)-\zeta\) and use continuity of translations. For \(p=\infty\), use uniform continuity of \(\zeta\). Choose \(J\) so that the norm of these principal blocks is at least \(c_0>0\). Then choose \(L\) sufficiently large that
\[
C_M\sum_{\ell\in\mathbb Z\setminus\{0\}}
2^{-(M-|s|)L|\ell|}<\frac{c_0}{2}.
\]
The reverse triangle inequality in each block \(j_m\) gives \(2^{j_ms}\|\Delta_{j_m}u_N\|_{L^p(\mathbb R^n)}\geq \displaystyle\frac{c_0}{2}\). The upper bound for all blocks follows by summing \eqref{eq:bloques-modulacion-compacta} and applying Young's inequality. Thus,
\[
cN^{\frac{1}{q}}\leq\|u_N\|_{B^s_{p,q}(\mathbb R^n)}\leq CN^{\frac{1}{q}}.
\]
By the choice of supports, the nonzero localizations of \(u_N\) are exactly \(v_1,\ldots,v_N\), and each has Besov norm between two positive constants. Consequently,
\[
cN^{\frac{1}{p}}\leq
\left\|\bigl(\|\eta_\nu u_N\|_{B^s_{p,q}(\mathbb R^n)}\bigr)_{\nu\in I}\right\|_{\ell^p(I)}
\leq CN^{\frac{1}{p}}.
\]
We interpret \(N^{\displaystyle\frac{1}{\infty}}=1\). Equivalence of the two norms with constants independent of \(N\) would require \(p=q\). For \(p=q<\infty\), the identity \(B^s_{p,p}(\mathbb R^n)=F^s_{p,p}(\mathbb R^n)\) and (a) give the asserted equivalence.
\end{proof}

\begin{remark}[Besov localization on manifolds]
Part (b) does not prevent localization of a Besov function with finitely many charts, but it does prevent using an unchanged $\ell^p(I)$ sum over an infinite atlas. On a noncompact manifold without boundary and of bounded geometry, one first defines the scale $F$, which does localize, and then
\[
B^s_{p,q}(M)
:=
(F^{s_0}_{p,p}(M),F^{s_1}_{p,p}(M))_{\theta,q},
\qquad
s=(1-\theta)s_0+\theta s_1.
\]
In the Euclidean case, Corollary~\ref{cor:interpolacion-real-BF-indices-finos} justifies this formula for any \(s_0<s<s_1\), with \(\theta=\displaystyle\frac{s-s_0}{s_1-s_0}\). In the construction on the manifold, independence of auxiliary parameters follows by combining localization formulas with reiteration. Uniform bounds on the charts and their transitions are the assumptions controlling the constants over an infinite atlas. The procedure will be revisited in Chapter~\ref{cap:geometria-acotada} and agrees with the Große--Schneider construction~\cite{traces} and Amann's localization principle~\cite{Amann2025FunctionSpaces}.
\end{remark}

\section{Open subsets, extension operators, and comparison of scales}
\label{sec:extension-dominios-fraccionarios}

For $A\in\{B,F\}$, define the restriction space
\[
A^s_{p,q}(\Omega)
:=
\{U\restriction_\Omega\mid U\in A^s_{p,q}(\mathbb R^n)\},
\qquad
\|u\|_{A^s_{p,q}(\Omega)}
:=
\inf_{\substack{U\in A^s_{p,q}(\mathbb R^n)\\U\restriction_\Omega=u}}\|U\|_{A^s_{p,q}(\mathbb R^n)}.
\]
This convention agrees with Definition~\ref{def:regularidad-intermedia-norma-r-n-abierto-s} when $A^s_{p,q}=F^s_{p,2}(\mathbb R^n)=H^{s,p}(\mathbb R^n)$.

Extension must use only values of the distribution inside the domain. To achieve this, we construct kernels supported in a cone that remains within the domain when subtracted from a point of the domain. Cancellation of all moments will allow the same kernels to be used at every regularity order.

\begin{lemma}[Conical reproducing kernels]
\label{lem:nucleos-rychkov-cono}
Let \(\omega\colon\mathbb R^{n-1}\longrightarrow\mathbb R\) be a Lipschitz function with constant \(L\), and set
\[
\Omega_\omega:=\{(x',x_n)\in\mathbb R^n\mid x_n>\omega(x')\},
\qquad
\mathcal K_L:=\{(y',y_n)\in\mathbb R^n\mid y_n<-(L+1)\|y'\|\}.
\]
There exist \(\varphi_0,\varphi,\psi_0,\psi\in\mathcal S(\mathbb R^n)\) such that, upon defining
\[
\varphi_j(x):=2^{jn}\varphi(2^jx),\qquad
\psi_j(x):=2^{jn}\psi(2^jx),\qquad j\in\mathbb N,
\]
the following properties hold.
\begin{enumerate}[label=(\alph*)]
\item The supports of \(\varphi_j,\psi_j\) are contained in \(\mathcal K_L\cap\{x_n\leq-c2^{-j}\}\), with \(c>0\). Moreover,
\[
\int_{\mathbb R^n}\varphi_0=\int_{\mathbb R^n}\psi_0=1,\qquad
\int_{\mathbb R^n}x^\alpha\varphi(x)\,dx
=\int_{\mathbb R^n}x^\alpha\psi(x)\,dx=0
\quad(\alpha\in\mathbb N_0^n).
\]
\item Denote by
\[
\mathcal S'(\Omega_\omega):=
\{U\restriction_{\Omega_\omega}\mid U\in\mathcal S'(\mathbb R^n)\}
\]
the space of restrictions of tempered distributions. For \(u\) in this space, \(\varphi_j*u\) is a well-defined smooth function on \(\Omega_\omega\), independent of the extension, and
\begin{equation}
u=\sum_{j\in\mathbb N_0}\psi_j*(\varphi_j*u)
\quad\text{in }\mathcal D'(\Omega_\omega).
\label{eq:reproduccion-conica-rychkov}
\end{equation}
\item For \(U\in\mathcal S'(\mathbb R^n)\), write
\[
Q_{j,\lambda}U(x):=
\sup_{y\in\mathbb R^n}
\frac{|(\varphi_j*U)(y)|}{(1+2^j\|x-y\|)^\lambda}.
\]
If \(\lambda>\displaystyle\frac{n}{a}\), where \(0<a<1\), then
\begin{equation}
\left\|\bigl(2^{js}\|Q_{j,\lambda}U\|_{L^p(\mathbb R^n)}\bigr)_{j\in\mathbb N_0}
\right\|_{\ell^q(\mathbb N_0)}
\leq C\|U\|_{B^s_{p,q}(\mathbb R^n)},
\label{eq:analisis-medias-locales-conicas-B}
\end{equation}
\begin{equation}
\left\|\bigl(2^{js}Q_{j,\lambda}U\bigr)_{j\in\mathbb N_0}
\right\|_{L^p(\mathbb R^n;\ell^q(\mathbb N_0))}
\leq C\|U\|_{F^s_{p,q}(\mathbb R^n)}
\label{eq:analisis-medias-locales-conicas}
\end{equation}
for the Banach indices defined in this chapter. The first bound also allows \(p=\infty\). If \(f_j=(\varphi_j*U)\restriction_{\Omega_\omega}\), then for every \(N>0\),
\begin{equation}
|\Delta_k[\psi_j*(\mathbf{1}_{\Omega_\omega}f_j)](x)|
\leq C_N2^{-N|j-k|}Q_{j,\lambda}U(x),
\qquad j,k\in\mathbb N_0.
\label{eq:casi-ortogonalidad-conica-rychkov}
\end{equation}
\end{enumerate}
\end{lemma}

\begin{proof}
First construct a function \(g\in\mathcal S(\mathbb R)\) supported in \([1,\infty)\), with integral one and all positive-order moments zero. Choose \(\rho\in C_c^\infty((1,2))\), with \(\displaystyle \int_{\mathbb R}\rho=1\). For \(0\leq\ell\leq N\), the polynomial interpolation formula at the points \(2^0,\ldots,2^N\), evaluated at zero, gives the coefficients
\[
a_{\ell,N}:=
\prod_{\substack{r\in\{0,\ldots,N\}\\r\neq\ell}}
\frac{2^r}{2^r-2^\ell}.
\]
Indeed, the polynomials
\[
\prod_{\substack{r\in\{0,\ldots,N\}\\r\neq\ell}}
\frac{t-2^r}{2^\ell-2^r}
\]
equal one at \(t=2^\ell\) and zero at the other points. Their combination with the values of a polynomial of degree at most \(N\) agrees with that polynomial, since their difference has \(N+1\) roots. In particular,
\[
\sum_{\ell=0}^{N}a_{\ell,N}2^{\ell m}
=
\begin{cases}1,&m=0,\\0,&1\leq m\leq N.\end{cases}
\]
Separating factors with \(r<\ell\) and \(r>\ell\) shows that
\[
|a_{\ell,N}|\leq C2^{-\frac{\ell(\ell+1)}{2}}.
\]
The constant is uniform because \(\displaystyle \prod_{r\in\mathbb N}(1-2^{-r})>0\). This last assertion follows by taking logarithms and using \(-\log(1-t)\leq2t\) for \(0\leq t\leq\displaystyle\frac{1}{2}\). For each \(\ell\), the products have a limit \(a_\ell\) as \(N\to\infty\). The preceding bound allows passage to the limit in every moment sum, even after multiplying by \(2^{\ell m}\). Thus,
\[
\sum_{\ell\in\mathbb N_0}a_\ell2^{\ell m}=\delta_{m0},
\qquad m\in\mathbb N_0.
\]
Define
\[
g(t):=\sum_{\ell\in\mathbb N_0}
a_\ell2^{-\ell}\rho(2^{-\ell}t).
\]
The seminorm \(\displaystyle \sup_{t\in\mathbb R}(1+|t|)^A|D^Bg(t)|\), with \(A,B\in\mathbb N_0\), is dominated by a constant times \(\displaystyle \sum_{\ell\in\mathbb N_0}|a_\ell|2^{\ell(A-B-1)}\). All these series converge. Thus \(g\) is Schwartz, has the required support, and its moments can be computed term by term.

Choose an invertible matrix \(A\) whose columns lie in a narrow cone around the positive vertical half-axis and whose positive cone lies in \(\{y_n>(L+1)\|y'\|\}\). This choice is obtained by taking \(n\) linearly independent vectors sufficiently close to the vertical vector. Set \(B:=-A\) and
\[
\vartheta(x):=
|\det B|^{-1}\prod_{i=1}^{n}g((B^{-1}x)_i),
\qquad
\vartheta_j(x):=2^{jn}\vartheta(2^jx),\quad j\in\mathbb Z.
\]
The change of variables \(x=By\) proves that \(\vartheta\) has integral one. Every positive-degree monomial in \(By\) expands as a sum of monomials in \(y\), each containing a positive-order moment of some factor \(g\). Thus all positive-order moments of \(\vartheta\) vanish. Its support lies in \(\mathcal K_L\), separated from the vertex by a plane \(x_n=-c\).

Now define
\[
\varphi_0:=\vartheta,\qquad
\psi_0:=3\vartheta-2\vartheta*\vartheta,
\]
and, for \(j\in\mathbb N\),
\[
\varphi_j:=\vartheta_j-\vartheta_{j-1},
\]
\[
\psi_j:=
3(\vartheta_j+\vartheta_{j-1})
-2(\vartheta_j*\vartheta_j+
\vartheta_j*\vartheta_{j-1}+
\vartheta_{j-1}*\vartheta_{j-1}).
\]
These are dilations of two fixed functions. Stability of the cone under addition preserves the supports. The convolution moment formula, obtained by expanding \((x+y)^\alpha\), shows that \(\varphi_j\) and \(\psi_j\) have all moments zero. In particular, the integral of \(\psi_j\) is \(6-6=0\).

The polynomial identity
\[
(a-b)\bigl(3(a+b)-2(a^2+ab+b^2)\bigr)
=(3a^2-2a^3)-(3b^2-2b^3)
\]
yields, in the convolution algebra, the telescoping sum
\[
\sum_{j=0}^{N}\psi_j*\varphi_j
=3\vartheta_N*\vartheta_N
-2\vartheta_N*\vartheta_N*\vartheta_N.
\]
Since \(\displaystyle \int_{\mathbb R^n}\vartheta=1\), \(\vartheta_N*\zeta\to\zeta\) in \(\mathcal S\) for \(\zeta\in\mathcal S\). This is checked by writing the difference as
\[
\int_{\mathbb R^n}\vartheta(y)
\bigl(\zeta(\,\cdot-2^{-N}y)-\zeta\bigr)\,dy
\]
and using continuity of translations in each Schwartz seminorm, dominated by a power of \(1+\|y\|\). By duality, the telescoping sum converges to the identity in \(\mathcal S'\). With the unitary Fourier normalization, this is equivalent to
\begin{equation}
\sum_{j\in\mathbb N_0}\widehat\psi_j(\xi)\widehat\varphi_j(\xi)
=(2\pi)^{-n},\qquad \xi\in\mathbb R^n.
\label{eq:bezout-diadica-conica}
\end{equation}

If \(x\in\Omega_\omega\) and \(y\in\mathcal K_L\), then
\[
x_n-y_n>\omega(x')+(L+1)\|y'\|
\geq\omega(x'-y')+\|y'\|.
\]
Thus \(x-y\in\Omega_\omega\). For a tempered extension \(U\) of \(u\), define \((\varphi_j*u)(x):=\langle U,\varphi_j(x-\,\cdot)\rangle\). The Schwartz function in the pairing has support inside \(\Omega_\omega\) at a positive distance from its boundary. Approximate it in \(\mathcal S\) by compact cutoffs whose supports remain in that domain. Two extensions agreeing on \(\Omega_\omega\) act identically on each cutoff and, after passage to the limit, on the whole function. This proves independence. Derivatives with respect to \(x\) are taken in the kernel, so the convolution is smooth. The same argument applies to \(\psi_j\), and restriction of the global identity already proved gives \eqref{eq:reproduccion-conica-rychkov}.

We prove the bounds. If \(\kappa_k\) is the kernel of \(\Delta_k\), cancellation and Schwartz decay imply, for any \(N,\lambda>0\),
\begin{equation}
\int_{\mathbb R^n}|(\kappa_k*\psi_j)(y)|
(1+2^j\|y\|)^\lambda\,dy
\leq C_{N,\lambda}2^{-N|j-k|}.
\label{eq:integral-ponderada-cancelacion}
\end{equation}
We give the scale estimate explicitly. Write
$a_j(x)=2^{jn}a(2^jx)$ and $b_k(x)=2^{kn}b(2^kx)$, with
$a,b\in\mathcal S$. If $j>k$ and all moments of $a$ vanish,
Taylor's formula with integral remainder gives
\[
|a_j*b_k(x)|\leq C_M2^{-jM}
\sum_{|\alpha|=M}\int_{\mathbb R^n}|a(z)|\,\|z\|^M
\int_0^1|D^\alpha b_k(x-t2^{-j}z)|\,dt\,dz.
\]
After $v=2^k x-t2^{k-j}z$, use
\[
1+2^j\|x\|
\leq 2^{j-k}(1+\|v\|)(1+\|z\|).
\]
Multiply by the weight, integrate, and change variables to obtain
\[
\begin{aligned}
\int_{\mathbb R^n}|a_j*b_k(x)|(1+2^j\|x\|)^\lambda\,dx
&\leq C_M2^{-(M-\lambda)(j-k)}
\left(\int |a(z)|\,\|z\|^M(1+\|z\|)^\lambda\,dz\right)\\
&\quad\cdot\sum_{|\alpha|=M}
\int |D^\alpha b(v)|(1+\|v\|)^\lambda\,dv.
\end{aligned}
\]
All integrals are finite. If $k>j$, interchange the kernels
and use the moments of $b$; the weight is at the coarse scale
and the same calculation gives $C_M2^{-M(k-j)}$, without
the cost $2^{\lambda(k-j)}$.
Taking $a=\psi$, $b$ equal to the kernel of the positive blocks,
and $M>N+\lambda$ proves the required bound for positive indices.
If one index is zero, use only cancellation at the positive
index; if both are zero, weighted integrability of two Schwartz
functions suffices.

The same estimate holds for \(\varphi_j*\widetilde\kappa_r\). Inserting the reconstruction of \(U\) and comparing weights at scales \(j\) and \(r\), we obtain, for any \(N>0\),
\begin{equation}
Q_{j,\lambda}U(x)
\leq C_N\sum_{r\in\mathbb N_0}
2^{-N|j-r|}P_{r,\lambda}(\Delta_rU)(x).
\label{eq:peetre-medias-locales}
\end{equation}
In detail, bound \(|\Delta_rU(z)|\) by \(P_{r,\lambda}(\Delta_rU)(x)(1+2^r\|x-z\|)^\lambda\). After dividing by the weight at scale \(j\), the triangle inequality for the weights introduces at most \(2^{\lambda|j-r|}\). Absorb this factor by taking a higher cancellation order in \eqref{eq:integral-ponderada-cancelacion}. The pairings and sums first converge in \(\mathcal S'\); the estimate holds with possibly infinite sides and, when the norm of \(U\) is finite, controls their pointwise convergence.

Take \(N>|s|\). Young's inequality for sequences and control of Peetre's maximal function by \(\mathcal M_a\) give \eqref{eq:analisis-medias-locales-conicas-B} and \eqref{eq:analisis-medias-locales-conicas}. For Besov spaces, use the scalar maximal operator on each block; for Triebel--Lizorkin spaces, use the exponents \(\displaystyle\frac{p}{a}\) and \(\displaystyle\frac{q}{a}\), both greater than one. The case \(q=\infty\) follows using the maximal operator of the supremum.

Finally,
\[
|\varphi_j*U(y)|
\leq Q_{j,\lambda}U(x)(1+2^j\|x-y\|)^\lambda.
\]
Inserting this bound into convolution with \(\kappa_k*\psi_j\) and using \eqref{eq:integral-ponderada-cancelacion} proves \eqref{eq:casi-ortogonalidad-conica-rychkov}. The indicator of the domain only decreases the integral of absolute values.
\end{proof}

This construction realizes Rychkov's conical reproduction \cite[Theorem~4.1]{Rychkov1999}. The distinction between \(\mathcal S'(\Omega)\) and \(\mathcal D'(\Omega)\), as well as the double-convolution formulation, can also be found in \cite[Section~1.1]{ShiYao2024}.

The preceding computation gives a property to be used when differentiating the extension. Let \(\eta_j,\theta_j\) and \(j\in\mathbb N\) be dyadic dilations of Schwartz functions with all moments zero. For a measurable set \(D\subseteq\mathbb R^n\) and \(m\geq0\), define
\[
T_mf:=
\sum_{j\in\mathbb N}2^{-jm}
\eta_j*\bigl(\mathbf{1}_D(\theta_j*f)\bigr).
\]
Estimates \eqref{eq:integral-ponderada-cancelacion} and \eqref{eq:peetre-medias-locales}, applied to these two families, prove
\begin{equation}
\|T_mf\|_{A^{s+m}_{p,q}}\leq C\|f\|_{A^s_{p,q}},
\qquad A\in\{B,F\}.
\label{eq:doble-convolucion-ganancia}
\end{equation}
Indeed, after multiplying the output block by \(2^{k(s+m)}\), the factor \(2^{-jm}\) leaves
\[
2^{k(s+m)}2^{-jm}2^{-N|j-k|}
\leq 2^{js}2^{-(N-|s+m|)|j-k|}.
\]
Choosing \(N>|s+m|\) allows the geometric summation in the proof to be repeated. The series converges in \(\mathcal S'\): a tempered distribution gives \(|\theta_j*f(x)|\leq C2^{jM_0}(1+\|x\|)^{M_0}\) for some \(M_0\), whereas pairing with \(\eta_j(-\,\cdot)*\zeta\) decays faster than every power of \(2^{-j}\). This last bound is exactly the preceding Taylor computation applied in Schwartz seminorms. The same justification applies to the series defining the extension operator.

\begin{theorem}[Extension on Lipschitz domains]
\label{teo:extension-universal-lipschitz-fraccionaria}
Let \(\Omega\subseteq\mathbb R^n\) be a bounded Lipschitz domain. There is a single linear operator \(E\), independent of the parameters, with \(Eu\restriction_\Omega=u\), inducing continuous operators
\[
E\colon A^s_{p,q}(\Omega)\longrightarrow A^s_{p,q}(\mathbb R^n),
\qquad A\in\{B,F\},\quad s\in\mathbb R,
\]
for all indices defined in this chapter. Moreover,
\[
E\colon W^{s,p}(\Omega)\longrightarrow W^{s,p}(\mathbb R^n)
\quad\text{is continuous if }s\geq0,\quad1<p<\infty.
\]
If \(s>0\) is not an integer, then
\[
W^{s,p}(\Omega)=B^s_{p,p}(\Omega)=F^s_{p,p}(\Omega)
\]
with equivalent norms. The constants may depend on the domain and the parameters, but the operator is the same.
\end{theorem}

\begin{proof}
Begin on \(\Omega_\omega\). For \(u\in\mathcal S'(\Omega_\omega)\), define
\[
E_\omega u:=
\sum_{j\in\mathbb N_0}
\psi_j*\bigl(\mathbf{1}_{\Omega_\omega}(\varphi_j*u)\bigr).
\]
The preceding lemma shows that this formula depends only on \(u\). For a tempered extension \(U\), the functions \(\varphi_j*U\) are smooth and of polynomial growth, so their product with the indicator defines a tempered distribution. The justification following \eqref{eq:doble-convolucion-ganancia} proves convergence of the series. If the evaluation point lies in \(\Omega_\omega\), the support of \(\psi_j\) allows the indicator to be removed. Conical reproduction then gives \(E_\omega u=u\) on that domain.

Multiply \eqref{eq:casi-ortogonalidad-conica-rychkov} by \(2^{ks}\), take \(N>|s|\), and apply Young's inequality and the local-means bounds. For any extension \(U\) in the scale under consideration, we obtain
\[
\|E_\omega u\|_{A^s_{p,q}(\mathbb R^n)}
\leq C\|U\|_{A^s_{p,q}(\mathbb R^n)}.
\]
Taking the infimum over \(U\restriction_{\Omega_\omega}=u\) proves the bound by the quotient norm.

We now relate this norm to the intrinsic norms. First let \(0<\sigma<1\) and \(u\in W^{\sigma,p}(\Omega_\omega)\). Its extension by zero belongs to \(L^p(\mathbb R^n)\), so \(u\in\mathcal S'(\Omega_\omega)\). For \(j\geq1\), cancellation of \(\varphi_j\) gives
\[
\varphi_j*u(x)=
\int_{\mathcal K_L}\varphi_j(y)\bigl(u(x-y)-u(x)\bigr)\,dy.
\]
Both arguments of \(u\) remain in the domain. Hölder's inequality with respect to the measure \(|\varphi_j(y)|\,dy\) and Tonelli's theorem imply
\[
\begin{aligned}
\sum_{j\in\mathbb N}2^{j\sigma p}
\|\varphi_j*u\|_{L^p(\Omega_\omega)}^p
&\leq C\int_{\mathcal K_L}
\left(\sum_{j\in\mathbb N}2^{j\sigma p}|\varphi_j(y)|\right)
\int_{\Omega_\omega}|u(x-y)-u(x)|^p\,dx\,dy\\
&\leq C[u]_{W^{\sigma,p}(\Omega_\omega)}^p.
\end{aligned}
\]
The second inequality uses
\[
\sum_{j\in\mathbb N}2^{j\sigma p}|\varphi_j(y)|
\leq C\|y\|^{-n-\sigma p}.
\]
This bound follows by separating \(2^j\|y\|\leq1\) and \(2^j\|y\|>1\) with a decay exponent greater than \(n+\sigma p\). The zero block is controlled by Young's inequality. For Besov synthesis, the unweighted version of \eqref{eq:integral-ponderada-cancelacion} and Young's inequality suffice, so
\[
\|E_\omega u\|_{B^\sigma_{p,p}(\mathbb R^n)}\leq
C\left(\|u\|_{L^p(\Omega_\omega)}
+[u]_{W^{\sigma,p}(\Omega_\omega)}\right).
\]
The global identification with \(W^{\sigma,p}(\mathbb R^n)\) proves extension in this first interval.

To treat higher-order derivatives, construct primitives preserving the cone. With the function \(g\) used above, set
\[
g_-(t):=\frac{1}{2}g(\frac{t}{2}),\qquad d(t):=g(t)-g_-(t).
\]
All moments of \(d\) vanish. For \(m\in\mathbb N\),
\[
h_m(t):=\frac{1}{(m-1)!}
\int_{-\infty}^{t}(t-r)^{m-1}d(r)\,dr
\]
satisfies \(D^mh_m=d\), is supported in \([1,\infty)\), and is Schwartz. Indeed, the vanishing moments allow the integral to be written as the negative of the integral from \(t\) to \(\infty\), proving decay of all its derivatives. Integrating by parts \(m\) times gives
\[
\int_{\mathbb R}t^\ell h_m(t)\,dt
=(-1)^m\frac{\ell!}{(\ell+m)!}
\int_{\mathbb R}t^{\ell+m}d(t)\,dt=0,
\qquad\ell\in\mathbb N_0.
\]

The telescoping identity for a product gives
\[
\prod_{r=1}^{n}g(y_r)-\prod_{r=1}^{n}g_-(y_r)
=\sum_{i=1}^{n}
\left(\prod_{r=1}^{i-1}g(y_r)\right)d(y_i)
\left(\prod_{r=i+1}^{n}g_-(y_r)\right),
\]
with the empty product equal to one. Replace \(d\) by \(D^mh_m\) and make the linear change \(x=By\). This gives functions \(\tau_{i,m}\in\mathcal S(\mathbb R^n)\) supported in the same cone, with all moments zero, such that
\[
\varphi=\sum_{i=1}^{n}D_{v_i}^{\,m}\tau_{i,m},
\qquad v_i:=Be_i.
\]
Here \(D_{v_i}\) is differentiation in the constant direction \(v_i\). If \(\tau_{i,m,j}(x):=2^{jn}\tau_{i,m}(2^jx)\), then
\[
\varphi_j*u=
2^{-jm}\sum_{i=1}^{n}
\tau_{i,m,j}*D_{v_i}^{\,m}u
\quad\text{in }\Omega_\omega.
\]
Integration by parts is distributional; the kernel supports remain inside the domain, so no boundary terms appear.

If \(u\in W^{m,p}(\Omega_\omega)\), extend each \(D_{v_i}^{\,m}u\) by zero to a function \(f_i\in L^p(\mathbb R^n)\). For \(|\alpha|\leq m\), differentiating the series for \(E_\omega u\) gives the zero block and the sums
\[
\sum_{j\in\mathbb N}2^{j(|\alpha|-m)}
\psi_j^\alpha*
\bigl(\mathbf{1}_{\Omega_\omega}(\tau_{i,m,j}*f_i)\bigr),
\qquad
\psi_j^\alpha:=2^{-j|\alpha|}D^\alpha\psi_j.
\]
Bound \eqref{eq:doble-convolucion-ganancia}, with scale \(F^0_{p,2}(\mathbb R^n)=L^p(\mathbb R^n)\), controls these sums in \(F^{m-|\alpha|}_{p,2}(\mathbb R^n)\hookrightarrow L^p(\mathbb R^n)\). The zero block is controlled by Young's inequality with \(D^\alpha\psi_0\). Since directional derivatives are finite combinations of partial derivatives of order \(m\), we conclude
\[
\|E_\omega u\|_{W^{m,p}(\mathbb R^n)}
\leq C\|u\|_{W^{m,p}(\Omega_\omega)}.
\]
The case \(m=0\) follows directly from the bound in \(F^0_{p,2}(\mathbb R^n)\), using the extension by zero of \(u\).

If \(s=m+\sigma\), \(0<\sigma<1\), each \(D_{v_i}^{\,m}u\) belongs to \(W^{\sigma,p}(\Omega_\omega)\). The first case of the proof gives an extension \(f_i\) in \(B^\sigma_{p,p}(\mathbb R^n)\) with controlled norm. The same representation, now without differentiating the output kernel, and \eqref{eq:doble-convolucion-ganancia} give
\[
\|E_\omega u\|_{B^{m+\sigma}_{p,p}(\mathbb R^n)}
\leq C\left(\|u\|_{L^p(\Omega_\omega)}
+\sum_{i=1}^{n}\|D_{v_i}^{\,m}u\|_{W^{\sigma,p}(\Omega_\omega)}\right)
\leq C\|u\|_{W^{m+\sigma,p}(\Omega_\omega)}.
\]
The global identification \(B^{m+\sigma}_{p,p}(\mathbb R^n)=W^{m+\sigma,p}(\mathbb R^n)\) completes the intrinsic bound at all noninteger orders.

For a bounded Lipschitz domain, choose a finite cover by boundary cylinders and an interior region. After rotations and translations, each boundary portion is the graph of a Lipschitz function. Extend that function to all of \(\mathbb R^{n-1}\) by \(\displaystyle \inf_{y\in D}(\omega(y)+L\|x'-y\|)\), where \(D\) is the graph domain; the Lipschitz inequality proves agreement with \(\omega\) on \(D\). This gives special domains \(\Omega_i\) agreeing with \(\Omega\) in smaller cylinders. Choose a smooth partition \(\eta_0,\ldots,\eta_N\) on a neighborhood of \(\overline\Omega\), with \(\eta_0\) supported in \(\Omega\), and cutoffs \(\chi_i\) equal to one near \(\operatorname{supp}\eta_i\) and supported in these smaller cylinders. Define
\[
E_\Omega u:=
(\eta_0u)^{\,0}
+\sum_{i=1}^{N}\chi_iE_i(\eta_i u).
\]
On each special domain, set \(\eta_i u\) equal to zero outside the cylinder; the cutoff vanishes near this artificial boundary. The superscript \(0\) indicates extension by zero of the interior part. On \(\Omega\), the formula equals \(\displaystyle \sum_{i=0}^{N}\eta_i u=u\). Smooth multipliers and finiteness of the cover give the bounds in restriction spaces.

For intrinsic norms, the Leibniz rule controls derivatives of the products. In the Gagliardo quotients, terms crossing an artificial boundary are estimated using the positive distance \(d\) from the cutoff support to that boundary and
\[
\int_{\{h\in\mathbb R^n\mid\|h\|\geq d\}}
\|h\|^{-n-\sigma p}\,dh=C_{n,\sigma,p}d^{-\sigma p}.
\]
Interior terms are controlled by the seminorm of \(u\) and by \(\displaystyle |\eta_i(x)-\eta_i(y)|\leq C\min\{1,\|x-y\|\}\). For higher-order products, the Leibniz rule also produces derivatives of \(u\) of order less than \(m\). Gluing at integer orders, which uses only the Leibniz rule, already gives a continuous operator from \(W^{1,p}(\Omega)\) to \(W^{1,p}(\mathbb R^n)\). Each of these lower derivatives belongs to \(W^{1,p}(\Omega)\); extend it with this operator and apply the increment lemma to control its \(W^{\sigma,p}(\Omega)\) norm. This justifies all fractional terms in the products and preserves the intrinsic bounds under gluing.

Finally, restricting a global function never increases either its derivative norms or its Gagliardo integral. This gives \(\|u\|_{W^{s,p}(\Omega)}\leq C\|u\|_{B^s_{p,p}(\Omega)}\). The constructed extension gives the reverse inequality. The identity \(B^s_{p,p}(\Omega)=F^s_{p,p}(\Omega)\) passes to the quotients and completes the proof.
\end{proof}

\begin{corollary}[Comparison of $H/W$ on a domain]
\label{cor:comparacion-H-W-dominios}
Let $\Omega$ be a bounded Lipschitz domain, let $s>0$ be noninteger, and let $1<p<\infty$. Then
\[
\begin{cases}
W^{s,p}(\Omega)\hookrightarrow H^{s,p}(\Omega),&1<p\leq2,\\
H^{s,p}(\Omega)\hookrightarrow W^{s,p}(\Omega),&2\leq p<\infty.
\end{cases}
\]
Equality holds for $p=2$ and, if $\Omega$ has nonempty interior, the embedding is strict for $p\neq2$. On an arbitrary open subset, only the embedding of the restriction space $B^s_{p,p}(\Omega)$ into the intrinsic Gagliardo space is always guaranteed; equality, even for $p=2$, requires an extension property.
\end{corollary}

\begin{proof}
Apply the extension operator from Theorem~\ref{teo:extension-universal-lipschitz-fraccionaria}, the global embedding from Corollary~\ref{cor:comparacion-H-W}, and restriction to the domain. The identity for \(p=2\) follows in the same way.

For strictness, fix a ball \(B\Subset\Omega\). The construction in Corollary~\ref{cor:comparacion-H-W} gives a function \(u\) supported in \(B\) that belongs to the larger scale but not the smaller one on \(\mathbb R^n\). Its restriction belongs to the larger scale on \(\Omega\). If it also belonged to the smaller scale, the extension operator would give an extension \(V\) in that global scale. Take \(\chi\in C_c^\infty(\Omega)\) equal to one near \(\overline B\). Since \(V=u\) on \(\Omega\), we have \(\chi V=u\) on \(\mathcal S'(\mathbb R^n)\). Continuity of multiplication by \(\chi\) would place \(u\) in the smaller global scale, a contradiction.

On an arbitrary open subset, restricting a global extension decreases the integrals defining the intrinsic norm. This gives \(B^s_{p,p}(\Omega)\hookrightarrow W^{s,p}(\Omega)\). Our proof of the reverse inclusion requires the continuous extension constructed for Lipschitz domains.
\end{proof}

\section{Sobolev embeddings of real order}
\label{sec:encajes-sobolev-orden-real}

The quantity
\[
\operatorname{ind}(s,p):=s-\frac{n}{p}
\]
is called the differential index. Bernstein's inequalities from Theorem~\ref{teo:desigualdades-bernstein} show that passage from $p_0$ to a larger exponent $p_1$ requires losing $n\left(\displaystyle\frac{1}{p_0}-\displaystyle\frac{1}{p_1}\right)$ derivatives. Interpolation explains the same rule starting from integer-order embeddings.

\begin{theorem}[Embeddings between scales]
\label{teo:encajes-escalas-orden-real}
Let $1<p_0\leq p_1<\infty$ and $s_0,s_1\in\mathbb R$.
\begin{enumerate}[label=(\alph*)]
\item If $s_0-\displaystyle\frac{n}{p_0}=s_1-\displaystyle\frac{n}{p_1}$, then
\[
B^{s_0}_{p_0,q}(\mathbb R^n)\hookrightarrow B^{s_1}_{p_1,q}(\mathbb R^n)
\]
for \(1\leq q\leq\infty\), and also
\[
F^{s_0}_{p_0,q}(\mathbb R^n)\hookrightarrow
F^{s_1}_{p_1,q}(\mathbb R^n)
\]
for the same fine indices.
\item If $s_0>s_1$ and $p$ is fixed, then
\[
A^{s_0}_{p,q_0}(\mathbb R^n)\hookrightarrow A^{s_1}_{p,q_1}(\mathbb R^n),
\qquad A\in\{B,F\},
\]
for any admissible fine indices $q_0,q_1$.
\item If $s_0=s_1$, the embedding within a single scale holds when $q_0\leq q_1$.
\end{enumerate}
\end{theorem}

\begin{proof}
For Besov spaces, Bernstein's inequality gives, for each \(j\in\mathbb N_0\),
\[
2^{js_1}\|\Delta_ju\|_{L^{p_1}(\mathbb R^n)}
\leq C2^{j(s_1+n(\frac{1}{p_0}-\frac{1}{p_1}))}
\|\Delta_ju\|_{L^{p_0}(\mathbb R^n)}
=C2^{js_0}\|\Delta_ju\|_{L^{p_0}(\mathbb R^n)}.
\]
Taking the \(\ell^q(\mathbb N_0)\) norm proves the first embedding.

For Triebel--Lizorkin spaces, the blocks must be estimated together before integration. If \(p_0=p_1\), then also \(s_0=s_1\), and the assertion is immediate. Assume \(p_0<p_1\) and set
\[
d=s_0-s_1=n(\frac{1}{p_0}-\frac{1}{p_1}),\qquad
f_j=2^{js_0}\Delta_ju,\qquad
G(x)=\|(f_j(x))_{j\in\mathbb N_0}\|_{\ell^q(\mathbb N_0)}.
\]
We have \(0<d<\displaystyle\frac{n}{p_0}\). A band-limited reproducing kernel gives \(f_j=\widetilde K_j*f_j\), where
\[
|\widetilde K_j(y)|\leq C_N2^{jn}
(1+2^j\|y\|)^{-N}.
\]
By Minkowski's inequality in \(\ell^q(\mathbb N_0)\), with the supremum if \(q=\infty\),
\[
\begin{split}
\|(2^{-jd}f_j(x))_{j\in\mathbb N_0}\|_{\ell^q(\mathbb N_0)}
&\leq C_N\int_{\mathbb R^n}
\sup_{j\in\mathbb N_0}
\bigl[2^{j(n-d)}(1+2^j\|y\|)^{-N}\bigr]G(x-y)\,dy\\
&\leq C\int_{\mathbb R^n}\|y\|^{d-n}G(x-y)\,dy
=:C I_dG(x).
\end{split}
\]
The second inequality follows because \(t^{n-d}(1+t)^{-N}\) is bounded for \(t>0\) if \(N>n-d\).

We prove the scalar bound needed for \(I_d\). For \(R>0\), split the integral into \(\|y\|<R\) and \(\|y\|\geq R\). Decomposing the first region into annuli \(2^{-k-1}R\leq\|y\|<2^{-k}R\), \(k\in\mathbb N_0\), and using the Hardy--Littlewood maximal operator gives
\[
\int_{\|y\|<R}\|y\|^{d-n}G(x-y)\,dy
\leq C R^d\mathcal MG(x)
\sum_{k\in\mathbb N_0}2^{-kd}
\leq C'R^d\mathcal MG(x).
\]
Hölder's inequality on the second region yields
\[
\int_{\|y\|\geq R}\|y\|^{d-n}G(x-y)\,dy
\leq C R^{d-\frac{n}{p_0}}\|G\|_{L^{p_0}(\mathbb R^n)},
\]
since \((n-d)p_0'>n\). When \(\mathcal MG(x)>0\), choose \(R=(\displaystyle\frac{\|G\|_{L^{p_0}(\mathbb R^n)}}{\mathcal MG(x)})^{\displaystyle\frac{p_0}{n}}\); if the maximal function vanishes, both integrals vanish. This gives
\[
I_dG(x)\leq C
\|G\|_{L^{p_0}(\mathbb R^n)}^{\frac{dp_0}{n}}
\bigl(\mathcal MG(x)\bigr)^{1-\frac{dp_0}{n}}.
\]
Since \(p_1(1-\displaystyle\frac{dp_0}{n})=p_0\), boundedness of the maximal operator on \(L^{p_0}(\mathbb R^n)\) implies \(\|I_dG\|_{L^{p_1}(\mathbb R^n)}\leq C\|G\|_{L^{p_0}(\mathbb R^n)}\). This inequality and the preceding vector-valued estimate prove \(F^{s_0}_{p_0,q}(\mathbb R^n)\hookrightarrow F^{s_1}_{p_1,q}(\mathbb R^n)\), including \(q=\infty\).

For (b), write \(\delta=s_0-s_1>0\). If \((b_j)_{j\in\mathbb N_0}\) is a nonnegative sequence,
\[
\|(2^{-j\delta}b_j)_{j\in\mathbb N_0}\|_{\ell^{q_1}(\mathbb N_0)}
\leq
\|(2^{-j\delta})_{j\in\mathbb N_0}\|_{\ell^{q_1}(\mathbb N_0)}
\sup_{j\in\mathbb N_0}b_j
\leq C_{\delta,q_1}\|(b_j)_{j\in\mathbb N_0}\|_{\ell^{q_0}(\mathbb N_0)}.
\]
Apply this inequality to \(b_j=2^{js_0}\|\Delta_ju\|_{L^p(\mathbb R^n)}\) for Besov spaces and pointwise to \(b_j(x)=2^{js_0}|\Delta_ju(x)|\) for Triebel--Lizorkin spaces. Finally, (c) is the inclusion \(\ell^{q_0}(\mathbb N_0)\hookrightarrow\ell^{q_1}(\mathbb N_0)\) when \(q_0\leq q_1\), applied in these same two positions.
\end{proof}

\begin{lemma}[Fractional Morrey estimate for the Gagliardo seminorm]
\label{teo:morrey-fraccionario-slobodeckij}
Let $0<s<1$, $1<p<\infty$, and $sp>n$. If $\alpha=s-\displaystyle\frac{n}{p}$, every $u\in W^{s,p}(\mathbb R^n)$ has a representative in $C_b^{0,\alpha}(\mathbb R^n)$, and
\[
\|u\|_{C_b^{0,\alpha}(\mathbb R^n)}
\leq
C\|u\|_{W^{s,p}(\mathbb R^n)}.
\]
\end{lemma}

\begin{proof}
Let $B=B(x,r)$. By Proposition~\ref{desigualdad de holder} and Theorem~\ref{teo:tonelli},
\begin{align*}
\int_B|u-u_B|^p
&\leq
\frac{1}{|B|}\int_B\int_B|u(y)-u(z)|^p\,dy\,dz\\
&\leq
C r^{-n}(2r)^{n+sp}[u]_{W^{s,p}(\mathbb R^n)}^p
=
C r^{n+p\alpha}[u]_{W^{s,p}(\mathbb R^n)}^p.
\end{align*}
For $B_k=B(x,2^{-k}r)$, Hölder's inequality and the preceding estimate give
\[
|u_{B_{k+1}}-u_{B_k}|
\leq |B_{k+1}|^{-\frac{1}{p}}\|u-u_{B_k}\|_{L^p(B_k)}
\leq C(2^{-k}r)^\alpha[u]_{W^{s,p}(\mathbb R^n)}.
\]
The geometric series yields a limit and the bound
\[
|\widetilde u(x)-u_{B(x,r)}|
\leq Cr^\alpha[u]_{W^{s,p}(\mathbb R^n)}.
\]
The limit is independent of the initial radius. Indeed, if
$2^{-k-1}r<\rho\leq2^{-k}r$, the two concentric balls have
comparable volumes, and Hölder's inequality on the larger one
gives $|u_{B(x,\rho)}-u_{B_k}|\leq C(2^{-k}r)^\alpha[u]$.
Thus all averages converge to the same limit as $\rho\to0^+$.
Theorem~\ref{teo:diferenciacion-lebesgue} identifies it with
$u(x)$ almost everywhere.

For $x\neq y$, take $r=\|x-y\|$ and $D=B(x,3r)$.
Both $B(x,r)$ and $B(y,r)$ lie in $D$ and have volumes
comparable with its volume. For either ball,
\[
|u_B-u_D|\leq |B|^{-\frac{1}{p}}\|u-u_D\|_{L^p(D)}
\leq Cr^\alpha[u]_{W^{s,p}(\mathbb R^n)}.
\]
The two limits of averages and these comparisons give
\[
|\widetilde u(x)-\widetilde u(y)|
\leq C\|x-y\|^\alpha[u]_{W^{s,p}(\mathbb R^n)}.
\]
Finally, $|u_{B(x,1)}|\leq C\|u\|_p$, uniformly in $x$.
The limit bound with $r=1$ controls the uniform norm of
$\widetilde u$ and proves the complete estimate.
\end{proof}

\begin{theorem}[Sobolev, Morrey, and the critical case]
\label{teo:encajes-sobolev-fraccionarios}
Let $s>0$ and $1<p<\infty$.
\begin{enumerate}[label=(\alph*)]
\item If $sp<n$ and $\displaystyle p_s^*=\displaystyle\frac{np}{n-sp}$, then
\[
H^{s,p}(\mathbb R^n)\hookrightarrow L^{p_s^*}(\mathbb R^n).
\]
If $s\notin\mathbb N$, also
\[
W^{s,p}(\mathbb R^n)
\hookrightarrow
L^{p_s^*,p}(\mathbb R^n)
\hookrightarrow
L^{p_s^*}(\mathbb R^n).
\]
\item If $sp=n$, both $H^{s,p}(\mathbb R^n)$ and $W^{s,p}(\mathbb R^n)$ embed into $L^q(\mathbb R^n)$ for every $p\leq q<\infty$, but not into $L^\infty(\mathbb R^n)$ in general. The limiting embedding
\[
B^{\frac{n}{p}}_{p,1}(\mathbb R^n)\hookrightarrow L^\infty(\mathbb R^n)
\]
is continuous.
\item If $sp>n$, the elements have continuous representatives. More precisely, if $k\in\mathbb N_0$ and $0<\alpha<1$ satisfy $k+\alpha<s-\displaystyle\frac{n}{p}$, then
\[
H^{s,p}(\mathbb R^n),\ W^{s,p}(\mathbb R^n)\hookrightarrow C_b^{k,\alpha}(\mathbb R^n).
\]
When $0<s<1$ and $sp>n$, the Slobodeckij scale reaches the exponent $\alpha=s-\displaystyle\frac{n}{p}$.
\end{enumerate}
The same results hold on a bounded Lipschitz domain after restriction; there, in the critical case, one obtains $L^q(\Omega)$ for every $1\leq q<\infty$.
\end{theorem}

\begin{proof}
We first explain the interpolation mechanism. By Theorem~\ref{teo:identificaciones-H-F-W-entero}, the embedding of the $F$ scale at constant differential index gives
\[
H^{s,p}(\mathbb R^n)=F^s_{p,2}(\mathbb R^n)\hookrightarrow F^0_{p_s^*,2}(\mathbb R^n)=L^{p_s^*}(\mathbb R^n).
\]
To obtain the Lorentz refinement without varying the spatial exponent in the real interpolation formula, choose
\[
0\leq s_0<s<s_1<\frac{n}{p},
\qquad
s=(1-\theta)s_0+\theta s_1,
\]
and define $q_i$ by $q_i^{-1}=p^{-1}-\displaystyle\frac{s_i}{n}$. The embedding just proved gives $H^{s_i,p}(\mathbb R^n)\hookrightarrow L^{q_i}(\mathbb R^n)$. Interpolate the identity operator using Theorem~\ref{teo:interpolacion-operadores-metodo-K} and use Theorems~\ref{teo:interpolacion-escalas-BF} and \ref{teo:interpolacion-Lp-lorentz}:
\[
W^{s,p}(\mathbb R^n)
=
(H^{s_0,p}(\mathbb R^n),H^{s_1,p}(\mathbb R^n))_{\theta,p}
\hookrightarrow
(L^{q_0}(\mathbb R^n),L^{q_1}(\mathbb R^n))_{\theta,p}
=
L^{p_s^*,p}(\mathbb R^n).
\]
Indeed, $(p_s^*)^{-1}=(1-\theta)q_0^{-1}+\theta q_1^{-1}$. Since $p<p_s^*$, monotonicity of the second index in Lorentz spaces, included in Theorem~\ref{teo:interpolacion-Lp-lorentz}, implies $L^{p_s^*,p}(\mathbb R^n)\hookrightarrow L^{p_s^*}(\mathbb R^n)$. This proves (a).

If $sp=n$, apply (a) with any $s-\varepsilon$ and choose $\varepsilon$ to obtain the desired finite exponent. On $\mathbb R^n$, this covers $q\geq p$; on a domain of finite measure, Proposition~\ref{desigualdad de holder} also covers $q<p$.

We show explicitly that the critical embedding into $L^\infty(\mathbb R^n)$ fails. Choose $\eta\in\mathcal S(\mathbb R^n)$ with $\eta(0)>0$ and
\[
\operatorname{supp}(\widehat\eta)
\subseteq
\left\{\xi\in\mathbb R^n\middle|\frac{3}{4}\leq\|\xi\|\leq\frac{4}{3}\right\},
\]
and define
\[
u_N(x):=\displaystyle\sum_{j=1}^N\frac{1}{j}\eta(2^jx).
\]
Finite interaction of the blocks and the change of variable $y=2^jx$ give
\[
\sup_{N\in\mathbb N}\|u_N\|_{B^{\frac{n}{p}}_{p,p}(\mathbb R^n)}
\leq
C\left(\displaystyle\sum_{j=1}^\infty\frac{1}{j^p}\right)^{\frac{1}{p}}<\infty.
\]
We develop the uniform bound for the critical Triebel--Lizorkin
norm. Set
\[
 a:=\frac{n}{p}>0,\qquad
 S_a:=\sum_{\ell=0}^{\infty}(\ell+2)2^{-a\ell}
 =\frac{2}{1-2^{-a}}+\frac{2^{-a}}{(1-2^{-a})^2},
 \qquad
 G_N(x):=\left(\sum_{j=1}^N\frac{2^{2ja}}{j^2}
                  |\eta(2^jx)|^2\right)^{\frac{1}{2}}.
\]
Fix $M>a$ and $C_M$ such that
$|\eta(z)|\leq C_M(1+\|z\|)^{-M}$. On the annulus
$2^{-k-1}<\|x\|\leq2^{-k}$, with $k\geq0$, use the fact that the
$\ell^2$ norm does not exceed the sum of absolute values. For
$1\leq j\leq k$, the inequality
$\displaystyle\frac{k+1}{j}\leq k-j+2$ gives
\[
 \sum_{j=1}^{\min\{k,N\}}\frac{2^{ja}}{j}|\eta(2^jx)|
 \leq\frac{\|\eta\|_{L^\infty}2^{ka}}{k+1}
        \sum_{\ell=0}^{\infty}(\ell+2)2^{-a\ell}
 =\frac{\|\eta\|_{L^\infty}S_a2^{ka}}{k+1}.
\]
If $k=0$, this sum is empty. For $j>k$,
$\|2^jx\|>2^{j-k-1}$ and $\displaystyle\frac{k+1}{j}\leq1$, so
\[
 \sum_{j=k+1}^N\frac{2^{ja}}{j}|\eta(2^jx)|
 \leq\frac{C_M2^M2^{ka}}{k+1}
         \sum_{\ell=1}^{\infty}2^{-(M-a)\ell}
 =\frac{C_M2^M2^{ka}}{k+1}
       \frac{2^{-(M-a)}}{1-2^{-(M-a)}}.
\]
An empty sum in this formula is interpreted as zero. Thus,
$G_N(x)\leq C\displaystyle\frac{2^{ka}}{k+1}$ on each annulus,
with $C$ independent of $k$ and $N$. On the ball $\|x\|\leq2^{-N}$,
the same calculation for the first sum, with $k=N$, gives
\[
 G_N(x)\leq\sum_{j=1}^N\frac{2^{ja}}{j}\|\eta\|_{L^\infty}
 \leq\frac{\|\eta\|_{L^\infty}S_a2^{Na}}{N+1}.
\]
For $\|x\|>1$, decay of $\eta$ implies
\[
 G_N(x)\leq C_M\|x\|^{-M}
             \sum_{j=1}^{\infty}2^{-j(M-a)}
 =D_M\|x\|^{-M},\qquad
 D_M:=C_M\frac{2^{-(M-a)}}{1-2^{-(M-a)}}.
\]
Let $\omega_n$ be the volume of the unit ball. Integrate over the
annuli with $0\leq k<N$, the ball of radius $2^{-N}$, and the exterior
of the unit ball. Since $ap=n$, $p>1$, and $Mp>n$, we obtain
\[
 \int_{\mathbb R^n}G_N(x)^p\,dx
 \leq \omega_n C^p\left(
          \sum_{k=0}^{N-1}\frac{1}{(k+1)^p}
          +\frac{1}{(N+1)^p}\right)
       +\frac{n\omega_nD_M^p}{Mp-n}
 \leq\omega_n C^p\left(2+\frac{1}{p-1}\right)
       +\frac{n\omega_nD_M^p}{Mp-n}.
\]
The final bound uses the comparison
$\sum_{j=1}^{\infty}j^{-p}\leq1+\int_1^\infty t^{-p}\,dt$.
Lemma~\ref{lem:sintesis-bloques-espectrales}, applied to the summands
with the stated spectral support, concludes
\[
 \sup_{N\in\mathbb N}\|u_N\|_{F^{\frac{n}{p}}_{p,2}(\mathbb R^n)}
 <\infty.
\]

However, $u_N(0)=\eta(0)\displaystyle\sum_{j=1}^N\displaystyle\frac{1}{j}\to\infty$. Theorem~\ref{teo:identificaciones-H-F-W-entero} always identifies $H^{\displaystyle\frac{n}{p},p}(\mathbb R^n)$ with $F^{\displaystyle\frac{n}{p}}_{p,2}(\mathbb R^n)$. If $\displaystyle\frac{n}{p}\notin\mathbb N$, Corollary~\ref{teo:identificaciones-H-W-B-F} identifies $W^{\displaystyle\frac{n}{p},p}(\mathbb R^n)$ with $B^{\displaystyle\frac{n}{p}}_{p,p}(\mathbb R^n)$; if $\displaystyle\frac{n}{p}\in\mathbb N$, the same Theorem~\ref{teo:identificaciones-H-F-W-entero} gives $W^{\displaystyle\frac{n}{p},p}(\mathbb R^n)=H^{\displaystyle\frac{n}{p},p}(\mathbb R^n)$. Therefore, in every case, neither $W^{\displaystyle\frac{n}{p},p}(\mathbb R^n)$ nor $H^{\displaystyle\frac{n}{p},p}(\mathbb R^n)$ embeds into $L^\infty(\mathbb R^n)$.

On the other hand, Theorem~\ref{teo:desigualdades-bernstein} and absolute summation give
\[
\|u\|_{L^\infty(\mathbb R^n)}
\leq
\displaystyle\sum_{j\in\mathbb N_0}\|\Delta_ju\|_{L^\infty(\mathbb R^n)}
\leq
C\displaystyle\sum_{j\in\mathbb N_0}2^{\frac{jn}{p}}\|\Delta_ju\|_{L^p(\mathbb R^n)}
=
C\|u\|_{B^{\frac{n}{p}}_{p,1}(\mathbb R^n)}.
\]

If $sp>n$, Theorem~\ref{teo:desigualdades-bernstein}, applied block by block, maps the source space into $B^{s-\displaystyle\frac{n}{p}}_{\infty,\infty}(\mathbb R^n)$, the Hölder--Zygmund scale. If $k+\alpha<s-\displaystyle\frac{n}{p}$, the uniform block sum and the estimate
\[
\|D^k\Delta_ju(\,\cdot+h)-D^k\Delta_ju\|_{L^\infty(\mathbb R^n)}
\leq
C\min\{1,2^j\|h\|\}\,2^{jk}\|\Delta_ju\|_{L^\infty(\mathbb R^n)}
\]
are split into $2^j\|h\|\leq1$ and $2^j\|h\|>1$; the two geometric series give the modulus $\|h\|^\alpha$. For $W^{s,p}(\mathbb R^n)$ with $0<s<1$, Lemma~\ref{teo:morrey-fraccionario-slobodeckij} directly provides the limiting exponent $s-\displaystyle\frac{n}{p}$. The extension from Theorem~\ref{teo:extension-universal-lipschitz-fraccionaria} transfers all embeddings to the domain.
\end{proof}

\begin{theorem}[Rellich--Kondrashov of real order]
\label{teo:rellich-orden-real}
Let $\Omega\subseteq\mathbb R^n$ be a bounded Lipschitz domain, and let $1<p<\infty$ and $s>t\geq0$. Then
\[
H^{s,p}(\Omega)\hookrightarrow H^{t,p}(\Omega)
\]
is compact. If $s,t$ are not integers, the same result holds for $W^{s,p}(\Omega)\hookrightarrow W^{t,p}(\Omega)$; with the usual integer-order interpretations, it holds for all nonnegative orders. Moreover, if $sp<n$ and $1\leq q<p_s^*$, the embeddings into $L^q(\Omega)$ are compact.
\end{theorem}

\begin{proof}
Let $(u_\nu)_{\nu\in\mathbb N}$ be a bounded sequence in one of the spaces of order $s$. Extend it by Theorem~\ref{teo:extension-universal-lipschitz-fraccionaria} and multiply by a cutoff $\zeta\in C_c^\infty(\mathbb R^n)$ equal to one on a neighborhood of $\overline\Omega$. This gives a sequence $(v_\nu)_{\nu\in\mathbb N}$ supported in a fixed compact set $K\subset\mathbb R^n$ and bounded in the corresponding global space. If
\[
P_J:=\displaystyle\sum_{j=0}^J\Delta_j,
\]
the block definition and Proposition~\ref{prop:young-sucesiones} give the uniform tail estimate
\[
\|(I-P_J)v_\nu\|_{A^{t}_{p,r}(\mathbb R^n)}
\leq
C2^{-J(s-t)}\|v_\nu\|_{A^{s}_{p,r}(\mathbb R^n)},
\qquad A\in\{B,F\}.
\]
Here take $A=F$, $r=2$ for the Bessel scale and $A=B$, $r=p$ for the noninteger Slobodeckij scale. We first prove compactness within each of these two scales, then treat integer or mixed orders of $W^{s,p}(\mathbb R^n)$.

Fix $J$. The operator $P_J$ is convolution with a kernel $K_J\in\mathcal S(\mathbb R^n)$. The sequence $(v_\nu)_{\nu\in\mathbb N}$ is bounded in $L^p(\mathbb R^n)$ and vanishes outside $K$. By Theorem~\ref{teo:young-convoluciones} and Proposition~\ref{prop:continuidad-traslaciones-Lp},
\[
\|P_Jv_\nu(\,\cdot+h)-P_Jv_\nu\|_{L^p(\mathbb R^n)}
\leq
\|K_J(\,\cdot+h)-K_J\|_{L^1(\mathbb R^n)}
\|v_\nu\|_{L^p(\mathbb R^n)}
\longrightarrow0
\]
uniformly in $\nu$ as $\|h\|\to0$. Moreover, Minkowski's integral inequality from Theorem~\ref{teo:minkowski-integral} gives
\[
\|P_Jv_\nu\|_{L^p(\{x\in\mathbb R^n\mid\|x\|>R\})}
\leq
\|v_\nu\|_{L^1(K)}
\sup_{y\in K}
\|K_J(\,\cdot-y)\|_{L^p(\{x\in\mathbb R^n\mid\|x\|>R\})}
\longrightarrow0
\]
uniformly in $\nu$ as $R\to\infty$; we also used $\|v_\nu\|_{L^1(K)}\leq |K|^{\displaystyle\frac{1}{p'}}\|v_\nu\|_{L^p(\mathbb R^n)}$. The Fréchet--Kolmogorov criterion from Proposition~\ref{corolario 14.47 relativamente compacto en Lp} shows that $(P_Jv_\nu)_{\nu\in\mathbb N}$ is relatively compact in $L^p(\mathbb R^n)$.

Since $P_Jv_\nu$ has frequencies in a fixed ball, Theorem~\ref{teo:desigualdades-bernstein} implies
\[
\|P_J(v_\nu-v_\mu)\|_{A^t_{p,r}(\mathbb R^n)}
\leq
C_J\|P_J(v_\nu-v_\mu)\|_{L^p(\mathbb R^n)}.
\]
Write $\displaystyle M_s:=\displaystyle\sup_{\nu\in\mathbb N}\|v_\nu\|_{A^s_{p,r}(\mathbb R^n)}<\infty$. For each $k\in\mathbb N$, choose $J_k$ so that
\[
2C2^{-J_k(s-t)}M_s<2^{-k}.
\]
Starting from the original sequence, recursively extract nested subsequences $(v_\nu^{(k)})_\nu$: once the subsequence at level $k-1$ is chosen, compactness of $(P_{J_k}v_\nu^{(k-1)})_\nu$ in $L^p(\mathbb R^n)$ and Bernstein's inequality allow us to choose $(v_\nu^{(k)})_\nu$ so that $(P_{J_k}v_\nu^{(k)})_\nu$ converges in $A^t_{p,r}(\mathbb R^n)$. Take the diagonal subsequence $w_k:=v_k^{(k)}$. For each fixed $k_0$, the tail $(w_k)_{k\geq k_0}$ is a subsequence of $(v_\nu^{(k_0)})_\nu$; thus $(P_{J_{k_0}}w_k)_{k\geq k_0}$ is Cauchy in the space of order $t$. Given $\varepsilon>0$, choose $k_0$ with $2^{-k_0}<\displaystyle\frac{\varepsilon}{2}$ and then $N\geq k_0$ so that for $k,\ell\geq N$,
\[
\|P_{J_{k_0}}(w_k-w_\ell)\|_{A^t_{p,r}(\mathbb R^n)}
<\frac{\varepsilon}{2}.
\]
The uniform estimate for the two tails then gives
\[
\|w_k-w_\ell\|_{A^t_{p,r}(\mathbb R^n)}<\varepsilon.
\]
Thus $(w_k)$ converges in the space of order $t$. Restriction to $\Omega$ is continuous, proving compactness between orders within each scale.

For integer or mixed orders of $W^{s,p}(\mathbb R^n)$, choose two nonintegers $a,b$ with $t<b<a<s$. The strict loss in order allows the fine block index to change. Indeed, if $A_j:=2^{js}\|\Delta_jv\|_{L^p(\mathbb R^n)}$, the source norm dominates $\displaystyle \sup_{j\in\mathbb N_0} A_j$, both for $B^s_{p,p}(\mathbb R^n)$ and for $F^s_{p,2}(\mathbb R^n)$. By convergence of the geometric series,
\[
\|v\|_{B^a_{p,p}(\mathbb R^n)}
\leq C\left(\sum_{j=0}^{\infty}
2^{-jp(s-a)}A_j^p\right)^{\frac{1}{p}}
\leq C_{s-a}\|v\|_{W^{s,p}(\mathbb R^n)}.
\]
When $s$ is an integer, use $W^{s,p}(\mathbb R^n)=H^{s,p}(\mathbb R^n)=F^s_{p,2}(\mathbb R^n)$ here, with equivalence of norms. At the target endpoint, the triangle inequality in $L^p(\mathbb R^n;\ell^2(\mathbb N_0))$ and Hölder's inequality for sequences give
\[
\|v\|_{F^t_{p,2}(\mathbb R^n)}
\leq C\sum_{j=0}^{\infty}2^{jt}\|\Delta_jv\|_{L^p(\mathbb R^n)}
\leq C_{b-t}\|v\|_{B^b_{p,p}(\mathbb R^n)}.
\]
If $t$ is not an integer, the inclusion $B^b_{p,p}(\mathbb R^n)\hookrightarrow B^t_{p,p}(\mathbb R^n)$ follows directly from the weights; if it is an integer, the last estimate and $F^t_{p,2}(\mathbb R^n)=W^{t,p}(\mathbb R^n)$ give the required inclusion. Applying these estimates to extensions with common support and using the compactness from $B^a_{p,p}(\mathbb R^n)$ to $B^b_{p,p}(\mathbb R^n)$ just proved gives a subsequence converging in $W^{t,p}(\Omega)$. This includes all nonnegative orders without identifying $B^j_{p,p}(\mathbb R^n)$ with $W^{j,p}(\mathbb R^n)$ at an integer order.

For the embedding into $L^q(\Omega)$, choose $0<t<s$ such that $q<p_t^*<p_s^*$ when $q>p$; if $q\leq p$, any $0<t<s$ and Proposition~\ref{desigualdad de holder} on the domain of finite measure suffice. In both cases, use
\[
A^{s}_{p,r}(\Omega)\hookrightarrow\hookrightarrow A^{t}_{p,r}(\Omega)\hookrightarrow L^q(\Omega).
\]
The second embedding is Theorem~\ref{teo:encajes-sobolev-fraccionarios}. This completes the proof.
\end{proof}

\section{Fractional-order traces}
\label{sec:trazas-fraccionarias-euclideas}

Restriction to the hyperplane \(\mathbb R^{n-1}\times\{0\}\) lowers dimension by one. A function concentrated in a strip of thickness \(2^{-j}\) has a spatial norm containing the factor \(2^{-\displaystyle\frac{j}{p}}\); evaluation at the center of the strip loses that factor. This explains the loss of \(\displaystyle\frac{1}{p}\) derivatives. In Triebel--Lizorkin spaces, integration in the normal direction also determines the fine index \(p\) of the target norm.

Assume \(n\geq2\). When \(n=1\), omit tangential variables and interpret the space on a point as \(\mathbb K\) with its absolute value; the same arguments give evaluation at the endpoints of an interval.

For a compact smooth hypersurface \(\Sigma\), define \(B^\sigma_{p,q}(\Sigma)\) using a finite atlas, a partition of unity, and the Euclidean norms of the local expressions. Lemma~\ref{lema: multiplicadores y difeomorfismos sobolev}, applied in dimension \(n-1\), compares any two such norms: insert one partition into the other and apply the coordinate-change operators on each overlap. Since there are finitely many terms, both inequalities have constants independent of the distribution. The resulting space is independent of the atlas.

\begin{lemma}[Dyadic trace estimate]
\label{lem:estimacion-diadica-traza}
Let \(1<p<\infty\) and \(j\in\mathbb N_0\), and let \(v\in L^p(\mathbb R^n)\) satisfy
\[
\operatorname{supp}\widehat v
\subseteq\{\xi\in\mathbb R^n\mid\|\xi\|\leq A2^j\}.
\]
Then
\begin{equation}
\|v(\,\cdot,0)\|_{L^p(\mathbb R^{n-1})}
\leq C_{A,p}2^{\frac{j}{p}}\|v\|_{L^p(\mathbb R^n)}.
\label{eq:traza-banda-limitada}
\end{equation}
Let \(s>\displaystyle\frac{1}{p}\), \(1\leq q\leq\infty\), \(u_j=\Delta_ju\), \(j\in\mathbb N_0\), and \(\sigma=s-\displaystyle\frac{1}{p}>0\). If the corresponding right-hand norm is finite, then
\begin{equation}
\left\|
\left(2^{\ell\sigma}
\left\|\Delta_\ell'\sum_{j\in\mathbb N_0}u_j(\,\cdot,0)\right\|_{L^p(\mathbb R^{n-1})}
\right)_{\ell\in\mathbb N_0}
\right\|_{\ell^q(\mathbb N_0)}
\leq C
\left\|\bigl(2^{js}\|u_j\|_{L^p(\mathbb R^n)}\bigr)_{j\in\mathbb N_0}\right\|_{\ell^q(\mathbb N_0)},
\label{eq:traza-diadica-B}
\end{equation}
and
\begin{equation}
\left(
\sum_{\ell\in\mathbb N_0}2^{\ell\sigma p}
\left\|\Delta_\ell'\sum_{j\in\mathbb N_0}u_j(\,\cdot,0)\right\|_{L^p(\mathbb R^{n-1})}^p
\right)^{\frac{1}{p}}
\leq C
\left\|\bigl(2^{js}u_j\bigr)_{j\in\mathbb N_0}\right\|_{L^p(\mathbb R^n;\ell^q(\mathbb N_0))}.
\label{eq:traza-diadica-F}
\end{equation}
Here \(\Delta_\ell'\) are blocks on \(\mathbb R^{n-1}\). Suprema replace the \(\ell^\infty(\mathbb N_0)\) norms. The constant in \eqref{eq:traza-diadica-F} can be chosen independently of \(q\in[1,\infty]\).
\end{lemma}

\begin{proof}
Choose a kernel \(\eta\in\mathcal S(\mathbb R)\) whose convolution multiplier \((2\pi)^{\displaystyle\frac{1}{2}}\widehat\eta\) equals one on \([-A,A]\). Reproduction in the normal variable gives
\[
v(x',0)=\int_{\mathbb R}2^j\eta(-2^jt)v(x',t)\,dt.
\]
Hölder's inequality in \(t\) implies
\[
|v(x',0)|^p
\leq C2^j\int_{\mathbb R}|v(x',t)|^p\,dt.
\]
Integration in \(x'\) proves \eqref{eq:traza-banda-limitada}. The formula is first obtained for Schwartz functions, then for \(v\) through reproduction by a kernel of enlarged frequency support. This kernel provides the smooth representative of \(v\) and allows passage to the limit in \(L^p(\mathbb R^n)\).

Set \(g_j=u_j(\,\cdot,0)\). Its tangential Fourier support lies in a ball of radius \(C2^j\), although it may contain low tangential frequencies. Thus,
\begin{equation}
\Delta_\ell'\sum_{j\in\mathbb N_0}g_j
=
\sum_{\substack{j\in\mathbb N_0\\j\geq\ell-N_0}}\Delta_\ell'g_j.
\label{eq:interacciones-traza}
\end{equation}
Young's inequality and \eqref{eq:traza-banda-limitada} give
\begin{equation}
2^{\ell\sigma}\|\Delta_\ell'g_j\|_{L^p(\mathbb R^{n-1})}
\leq C2^{-(j-\ell)\sigma}2^{js}\|u_j\|_{L^p(\mathbb R^n)},
\qquad j\geq\ell-N_0.
\label{eq:nucleo-hardy-traza}
\end{equation}
The sequence
\(\bigl(2^{-m\sigma}\mathbf{1}_{\{m\geq-N_0\}}\bigr)_{m\in\mathbb Z}\)
is summable. Young's inequality for sequences proves \eqref{eq:traza-diadica-B}.

For the second estimate, integrate each block over a different strip. Define
\[
I_j:=(2^{-j},2^{1-j}],\qquad j\in\mathbb N_0.
\]
These intervals are disjoint and have length \(2^{-j}\). If \(t\in I_j\), Peetre's definition implies
\[
|u_j(x',0)|\leq3^\lambda
P_{j,\lambda}u_j(x',t).
\]
Multiply by \(2^{jsp}\) and integrate in \(I_j\) and then in \(x'\). Summing in \(j\), disjointness of the intervals gives
\[
\begin{aligned}
\sum_{j\in\mathbb N_0}2^{j(sp-1)}\|g_j\|_{L^p(\mathbb R^{n-1})}^p
&\leq C
\sum_{j\in\mathbb N_0}\int_{\mathbb R^{n-1}\times I_j}
\bigl(2^{js}P_{j,\lambda}u_j(x)\bigr)^p\,dx\\
&\leq C\int_{\mathbb R^n}
\sup_{j\in\mathbb N_0}
\bigl(2^{js}P_{j,\lambda}u_j(x)\bigr)^p\,dx.
\end{aligned}
\]
Choose \(0<a<1\) and \(\lambda>\displaystyle\frac{n}{a}\). Peetre control and monotonicity of the maximal operator give
\[
\sup_{j\in\mathbb N_0}2^{js}P_{j,\lambda}u_j
\leq C\left[
\mathcal M\left(\sup_{j\in\mathbb N_0}|2^{js}u_j|^a\right)
\right]^{\frac{1}{a}}.
\]
Since \(\displaystyle\frac{p}{a}>1\), boundedness of the maximal operator implies
\[
\left(\sum_{j\in\mathbb N_0}2^{j\sigma p}\|g_j\|_{L^p(\mathbb R^{n-1})}^p\right)^{\frac{1}{p}}
\leq C\left\|\sup_{j\in\mathbb N_0}|2^{js}u_j|\right\|_{L^p(\mathbb R^n)}
\leq C\left\|\bigl(2^{js}u_j\bigr)_{j\in\mathbb N_0}\right\|_{L^p(\mathbb R^n;\ell^q(\mathbb N_0))}.
\]
Finally, \eqref{eq:interacciones-traza}, Young's convolution inequality, and the same geometric sequence from \eqref{eq:nucleo-hardy-traza} transfer this norm of the \(g_j\) to the \(B^\sigma_{p,p}(\mathbb R^{n-1})\) norm of their sum. This proves \eqref{eq:traza-diadica-F} with a constant independent of \(q\).

The trace series converge in \(L^p(\mathbb R^{n-1})\): in both situations, the sequence \((2^{j\sigma}\|g_j\|_{L^p(\mathbb R^{n-1})})_{j\in\mathbb N_0}\) is bounded and \(\displaystyle \sum_{j\in\mathbb N_0}2^{-j\sigma}<\infty\). Thus all steps for finite families pass to the limit.
\end{proof}

\begin{theorem}[Hyperplane trace and coretraction]
\label{teo:traza-fraccionaria-euclidea-hiperplano}
Let \(1<p<\infty\), \(1\leq q\leq\infty\), and \(s>\displaystyle\frac{1}{p}\). Restriction to the hyperplane induces continuous surjective operators
\begin{equation}
\operatorname{tr}\colon B^s_{p,q}(\mathbb R^n)
\longrightarrow B^{s-\frac{1}{p}}_{p,q}(\mathbb R^{n-1}),
\qquad
\operatorname{tr}\colon F^s_{p,q}(\mathbb R^n)
\longrightarrow B^{s-\frac{1}{p}}_{p,p}(\mathbb R^{n-1}).
\label{eq:trazas-BF-hiperplano}
\end{equation}
Both have a continuous linear coretraction. In particular,
\begin{equation}
\operatorname{tr}\colon H^{s,p}(\mathbb R^n)
\longrightarrow B^{s-\frac{1}{p}}_{p,p}(\mathbb R^{n-1})
\label{eq:traza-H-hiperplano}
\end{equation}
is a retraction.
\end{theorem}

\begin{proof}
For each element of the spaces in the statement, define
\[
\operatorname{tr}u:=\sum_{j\in\mathbb N_0}(\Delta_ju)(\,\cdot,0).
\]
The lemma proves convergence and continuity. This definition agrees with canonical evaluation in the normal direction. Indeed, \eqref{eq:traza-banda-limitada} holds on every hyperplane of height \(t\) with the same constant, and
\[
\sup_{t\in\mathbb R}\|(\Delta_ju)(\,\cdot,t)\|_{L^p(\mathbb R^{n-1})}
\leq C2^{-j(s-\frac{1}{p})}\|u\|_{A^s_{p,q}}.
\]
Each block is a continuous function of \(t\) with values in \(L^p(\mathbb R^{n-1})\). This follows from one-dimensional reproduction and continuity of translations of its kernel. The series converges uniformly in this space of continuous functions. Its sum represents \(u\) as a distribution by dyadic reconstruction, and its value at zero is the defined trace. Thus the trace is independent of the resolution. If \(u\) vanishes on the half-space \(t>0\), its continuous representative in \(t\) vanishes for \(t>0\) and also at \(t=0\). This observation will give independence from extensions.

Choose \(\rho\in\mathcal S(\mathbb R)\) with \(\rho(0)=1\) and compact Fourier support. For data \(g\) in the target space, set
\begin{equation}
(\operatorname{ex}g)(x',t):=
\sum_{j\in\mathbb N_0}\rho(2^jt)\Delta_j'g(x').
\label{eq:extension-diadica-traza}
\end{equation}
The summands have spectral support in annuli of \(\mathbb R^n\) comparable to \(2^j\), except for the zero block, whose support lies in a ball. Moreover,
\begin{equation}
\|\rho(2^j\,\cdot)\Delta_j'g\|_{L^p(\mathbb R^n)}
=2^{-\frac{j}{p}}\|\rho\|_{L^p(\mathbb R)}
\|\Delta_j'g\|_{L^p(\mathbb R^{n-1})}.
\label{eq:norma-bloque-extension}
\end{equation}
Lemma~\ref{lem:sintesis-bloques-espectrales} gives the Besov bound and distributional convergence.

For Triebel--Lizorkin spaces, we prove the stronger bound with fine index one; the inclusion \(\ell^1(\mathbb N_0)\hookrightarrow\ell^q(\mathbb N_0)\) will give the other indices. Write
\[
b_j(x'):=2^{j(s-\frac{1}{p})}|\Delta_j'g(x')|,\qquad j\in\mathbb N_0.
\]
Spectral synthesis and decay of \(\rho\) give
\[
\|\operatorname{ex}g\|_{F^s_{p,1}(\mathbb R^n)}
\leq C
\left\|\sum_{j\in\mathbb N_0}
2^{\frac{j}{p}}(1+2^j|t|)^{-M}b_j(x')\right\|_{L^p(\mathbb R^{n-1}\times\mathbb R,dx'\,dt)}.
\]
Take \(M>\displaystyle\frac{1}{p}\) and split \(\mathbb R\setminus\{0\}\) into
\[
A_m:=\{t\in\mathbb R\mid2^{-m-1}<|t|\leq2^{-m}\},
\qquad m\in\mathbb Z.
\]
Their measure is \(2^{-m}\). For \(t\in A_m\), multiplying the sum by \(2^{-\displaystyle\frac{m}{p}}\) bounds it by
\[
C\sum_{j\in\mathbb N_0}c_{j-m}b_j(x'),\qquad
c_\ell:=
\begin{cases}
2^{\frac{\ell}{p}},&\ell\leq0,\\
2^{-\ell(M-\frac{1}{p})},&\ell>0.
\end{cases}
\]
The sequence \((c_\ell)_{\ell\in\mathbb Z}\) belongs to \(\ell^1(\mathbb Z)\). Young's inequality in the indices \(j,m\), followed by Tonelli's theorem in \(x'\), yields
\begin{equation}
\|\operatorname{ex}g\|_{F^s_{p,q}(\mathbb R^n)}^p
\leq C\sum_{j\in\mathbb N_0}
2^{jp(s-\frac{1}{p})}\|\Delta_j'g\|_{L^p(\mathbb R^{n-1})}^p
=C\|g\|_{B^{s-\frac{1}{p}}_{p,p}(\mathbb R^{n-1})}^p.
\label{eq:cota-extension-F}
\end{equation}
The constants are independent of the number of summands.

Finally, \(s-\displaystyle\frac{1}{p}>0\) implies \(\displaystyle \sum_{j\in\mathbb N_0}\|\Delta_j'g\|_{L^p(\mathbb R^{n-1})}<\infty\). The series in \eqref{eq:extension-diadica-traza} converges uniformly in \(t\) with values in \(L^p(\mathbb R^{n-1})\) and equals \(\displaystyle \sum_{j\in\mathbb N_0}\Delta_j'g=g\) at \(t=0\). The canonical interpretation of the trace proves \(\operatorname{tr}\operatorname{ex}=I\).
\end{proof}

\begin{theorem}[Fractional trace on smooth domains]
\label{teo:traza-fraccionaria-euclidea}
Let \(\Omega\subseteq\mathbb R^n\) be a bounded domain with boundary \(C^\infty\), and let \(1<p<\infty\), \(1\leq q\leq\infty\), and \(s>\displaystyle\frac{1}{p}\). There are continuous surjective operators
\[
\operatorname{Tr}\colon B^s_{p,q}(\Omega)
\longrightarrow B^{s-\frac{1}{p}}_{p,q}(\partial\Omega),\qquad
\operatorname{Tr}\colon F^s_{p,q}(\Omega)
\longrightarrow B^{s-\frac{1}{p}}_{p,p}(\partial\Omega).
\]
Each has a continuous linear right inverse. In particular,
\[
\operatorname{Tr}\colon H^{s,p}(\Omega)
\longrightarrow B^{s-\frac{1}{p}}_{p,p}(\partial\Omega)
\]
is a retraction; the same holds with \(W^{s,p}(\Omega)\) as source space when \(s\) is not an integer.
\end{theorem}

\begin{proof}
Choose an extension \(U\) of \(u\) to \(\mathbb R^n\), boundary charts \(\Phi_i\colon Q\longrightarrow U_i\), \(i\in\{1,\ldots,N\}\), and subordinate smooth cutoffs \(\eta_i\) whose sum equals one near \(\partial\Omega\). The charts satisfy
\[
\Phi_i(Q\cap\{t>0\})=U_i\cap\Omega,\qquad
\Phi_i(Q\cap\{t=0\})=U_i\cap\partial\Omega.
\]
Choose the cutoffs compactly supported inside each chart. Define
\begin{equation}
\operatorname{Tr}_\Omega u:=
\sum_{i=1}^{N}
\left[\operatorname{tr}\bigl((\eta_iU)\circ\Phi_i\bigr)\right]
\circ(\Phi_i\restriction_{\{t=0\}})^{-1},
\label{eq:traza-dominio-local}
\end{equation}
with extension by zero of each localized expression. The multiplication and composition lemmas and the hyperplane theorem bound the sum by \(C\|U\|_{A^s_{p,q}}\).

If \(U_1,U_2\) extend the same \(u\), the localized expression of \(U_1-U_2\) vanishes on \(t>0\). The final observation about the canonical trace in the preceding theorem shows that its trace is zero. Thus the formula does not depend on the extension. Taking the infimum over all extensions gives the bound by the quotient norm of \(u\).

To compare two systems of charts, fix \(0<\varepsilon<s-\displaystyle\frac{1}{p}\). Corollary~\ref{cor:aproximacion-escalas-BF} allows approximation of \(U\), after localization, in \(B^{s-\varepsilon}_{p,p}(\mathbb R^n)\) by smooth functions. The two trace constructions agree for these functions. Multiplication, composition, and trace are continuous at this lower order, which still exceeds \(\displaystyle\frac{1}{p}\). Passing to the limit, the resulting traces agree as distributions on the boundary. Since both already belong to the target space of order \(s-\displaystyle\frac{1}{p}\), this identity is also equality of their elements in that space. The argument covers \(q=\infty\).

We now construct a right inverse. Choose cutoffs \(\zeta_i\in C^\infty(\partial\Omega)\) supported inside the charts, with
\(\displaystyle \sum_{i=1}^{N}\zeta_i^2=1\)
To obtain them, start with cutoffs whose sum is positive and divide each by the square root of the sum of their squares. Given \(g\), express
\[
g_i:=(\zeta_i g)\circ(\Phi_i\restriction_{\{t=0\}})
\]
on \(\mathbb R^{n-1}\), extending it by zero. Apply \(\operatorname{ex}\) from \eqref{eq:extension-diadica-traza}. Choose on \(Q\) a cutoff \(\chi_i\) whose restriction to \(t=0\) is the expression of \(\zeta_i\), and extend it smoothly with compact support inside \(Q\). The sum of the functions \(\chi_i\operatorname{ex}g_i\) transported by \(\Phi_i^{-1}\) and restricted to \(\Omega\) defines a continuous operator. Its trace is
\[
\sum_{i=1}^{N}\zeta_i(\zeta_i g)=g.
\]
The Bessel and Slobodeckij identifications give the final assertions.
\end{proof}

\begin{remark}[Lipschitz boundaries]
\label{obs:traza-lipschitz-alcance}
On a Lipschitz boundary, coordinate changes do not have the derivatives used in the preceding theorem. The trace formulation at this boundary regularity requires a separate development; in the interval
\[
\frac{1}{p}<s<1+\frac{1}{p}
\]
the target space has order between zero and one and is described through differences in boundary charts. For the hyperplane and smooth boundaries, see, respectively,
\cite[Sections~4.4.2 and~5.4]{Triebel2}. 
\end{remark}

\part{Nonlinear analysis on Riemannian manifolds}
\chapter*{Introduction to Part III}
\addcontentsline{toc}{chapter}{Introduction to Part III}
\markboth{Part III. Nonlinear analysis on Riemannian manifolds}{Introduction to Part III}

We now return to manifolds with the tools developed in the Euclidean setting. This is one of the central steps in the book. It is not enough to replace partial derivatives with covariant derivatives and repeat the same formulas: working on a manifold requires the Riemannian measure, bundle metrics, connections, supports of sections, and transformations between charts and trivializations. Each of these structures determines how an analytic notion is correctly formulated.

The first chapters construct Sobolev spaces of functions and of sections of vector bundles, and develop the integration-by-parts formulas used to define weak covariant derivatives and formal adjoints. Differential operators on bundles and vector-valued distributions are then introduced. The principal symbol and ellipticity allow regularity to be recovered from the equation. Later, the connection Laplacian and the heat semigroup provide intrinsic scales of fractional regularity.

The final stage of this development requires us to leave behind the convenience of compactness. On a noncompact manifold, a local estimate can be globalized only if its constants remain controlled as the chart moves. Bounded geometry provides precisely this uniform control. Normal and Fermi coordinates, partitions of unity, and uniform trivializations then allow the results to be extended to submanifolds and manifolds with boundary. Throughout these chapters, we will continually alternate between two operations: expressing a geometric object in coordinates in order to calculate, and reconstructing a global, intrinsic statement from local estimates.

\begin{semblanzaHistorica}{Aubin and analysis on manifolds}
The development of geometric equations in the twentieth century made it necessary to adapt the tools of functional analysis to spaces whose geometry could not be ignored. Thierry Aubin played a central role in systematizing Sobolev inequalities and variational problems on manifolds, particularly in the study of conformal equations. Subsequently, Hebey, Strichartz, Eichhorn, Grigor'yan, Amann, and many others extended these ideas to complete and noncompact manifolds, manifolds of bounded geometry, manifolds with boundary, and spaces with more general structures. In this part we follow the same course: beginning with the local theory and identifying the geometric hypotheses that turn it into a global theory.
\end{semblanzaHistorica}

\chapter{Sobolev spaces on Riemannian manifolds}

The Dirichlet energy provides an intrinsic starting point for transferring Sobolev spaces to a Riemannian manifold $(M,\mathbf{g})$. For scalar data $f$, we formally consider
\[
\mathcal E_{\mathbf{g}}(u)
=
\frac{1}{2}\int_M |du|_{\mathbf{g}}^2\,d\lambda_{\mathbf{g}}
-\int_M fu\,d\lambda_{\mathbf{g}}.
\]
For smooth $u$ and $\varphi$, with interior support or boundary conditions that eliminate the boundary term, the first variation gives
\[
\int_M \langle du,d\varphi\rangle_{\mathbf{g}}\,d\lambda_{\mathbf{g}}
=
\int_M f\varphi\,d\lambda_{\mathbf{g}}.
\]
This identity is the weak formulation of $-\Delta_{\mathbf{g}} u=f$ and remains meaningful when $u$ has only one square-integrable weak derivative.

The metric enters both the norm of $du$ and the measure $\lambda_{\mathbf{g}}$, but the definition must be independent of coordinates. Charts allow a function to be compared locally with its Euclidean representative; partitions of unity assemble these comparisons, and covariant derivatives provide a global description of higher orders. On a compact manifold, finiteness of the cover makes the norms obtained from these constructions equivalent.

This chapter will define the spaces $W^{m,p}(M)$, prove their completeness, and establish independence from the relevant auxiliary data. It will then transfer the density, embedding, and Rellich--Kondrashov theorems to the manifold. The following example first shows how the same geometric operator takes different expressions in different charts while still representing the Laplace--Beltrami operator.

\begin{example}[The Laplace--Beltrami operator on the sphere in different coordinates]

Let $S^{2}\subset\mathbb{R}^{3}$ be the unit sphere and let $\iota\colon S^{2}\hookrightarrow\mathbb{R}^{3}$ be the canonical inclusion. Denote the Euclidean metric on $\mathbb{R}^{3}$ by $\overline{\mathbf{g}}$. The Riemannian metric on $S^{2}$ is defined as the pullback of $\overline{\mathbf{g}}$ by the inclusion:
\[
\mathbf{g}=\iota^{*}\overline{\mathbf{g}},
\qquad
\mathbf{g}_{p}(X_{p},Y_{p})=\overline{\mathbf{g}}_{\iota(p)}(d\iota_{p}(X_{p}),d\iota_{p}(Y_{p})),
\quad X_{p},Y_{p}\in T_{p}S^{2}.
\]
If $(x^{1},x^{2})$ are local coordinates induced by a chart $\phi\colon U\subseteq S^{2}\longrightarrow V\subseteq\mathbb{R}^{2}$, denote by
$\Phi:=\iota\circ\phi^{-1}\colon V=\phi(U)\longrightarrow\mathbb{R}^{3}$
the composition of the inverse chart with the inclusion. By the definition of the induced metric,
\[
\mathbf{g}=\iota^{*}\overline{\mathbf{g}},
\qquad
\mathbf{g}_{p}(X_{p},Y_{p})
=\overline{\mathbf{g}}_{\iota(p)}\big(d\iota_{p}(X_{p}),  d\iota_{p}(Y_{p})\big),
\quad X_{p},Y_{p}\in T_{p}S^{2}.
\]
In particular, for the coordinate vectors
\[
\boldsymbol{\partial}_{i}\big|_{p}
=(d\phi^{-1})_{\phi(p)} \left(\boldsymbol{\partial}_{i}\bigg|_{\phi(p)}\right),
\]
we have
\[
d\iota_{p} \left(\boldsymbol{\partial}_{i}\big|_{p}\right)
=d\iota_{p}\circ \left((d\phi^{-1})_{\phi(p)}\left(\boldsymbol{\partial}_{i}\bigg|_{\phi(p)}\right)\right)
=d(\iota\circ\phi^{-1})_{\phi(p)} \left(\boldsymbol{\partial}_{i}\bigg|_{\phi(p)}\right)
=d\Phi_{\phi(p)} \left(\boldsymbol{\partial}_{i}\bigg|_{\phi(p)}\right).
\]
Thus, evaluating the definition of $\mathbf{g}$ at the point $p$,
\[
g_{ij}(p)
=\mathbf{g}_{p} \left(
\boldsymbol{\partial}_{i}\big|_{p},
\boldsymbol{\partial}_{j}\big|_{p}
\right)
=\overline{\mathbf{g}}_{\iota(p)} \left(
d\Phi_{\phi(p)} \left(\boldsymbol{\partial}_{i}\bigg|_{\phi(p)}\right),
d\Phi_{\phi(p)} \left(\boldsymbol{\partial}_{j}\bigg|_{\phi(p)}\right)
\right).
\]
Finally, since $d\Phi_{\phi(p)}\left(\displaystyle\boldsymbol{\partial}_{i}\Biggr|_{\phi(p)}\right)$ agrees with the partial derivative $\displaystyle\frac{\partial\Phi}{\partial x_{i}}(\phi(p))$ at $\mathbb{R}^{3}$, we obtain
\[
g_{ij}(p)
=\left\langle
\displaystyle\frac{\partial\Phi}{\partial x^{i}}(\phi(p)),
\displaystyle\frac{\partial\Phi}{\partial x^{j}}(\phi(p))
\right\rangle_{\mathbb{R}^{3}},
\]
where $\langle\cdot,\cdot\rangle_{\mathbb{R}^{3}}$ is the Euclidean inner product. This identity shows that the local coefficients of the induced metric on $S^{2}$ are computed, for each $p\in U$, as the inner products of the partial derivatives of the parametrization $\Phi=\iota\circ\phi^{-1}$ evaluated at $\phi(p)\in\mathbb{R}^{2}$.

We shall use spherical coordinates $(\theta,\varphi)$ with the following convention: $\theta$ is the radial angle in the $xy$ plane, $0<\theta<2\pi$, and $\varphi$ is the azimuthal angle measured from the positive $z$ axis, $0<\varphi<\pi$. The inverse of these coordinates defines a chart on an open subset of $S^{2}$. The corresponding map is
\[
\Phi(\theta,\varphi)=\bigl(\sin\varphi\cos\theta,\ \sin\varphi\sin\theta,\ \cos\varphi\bigr).
\]
Its partial derivatives are
\[
\partial_{\theta}\Phi=\bigl(-\sin\varphi\sin\theta,\ \sin\varphi\cos\theta,\ 0\bigr),
\qquad
\partial_{\varphi}\Phi=\bigl(\cos\varphi\cos\theta,\ \cos\varphi\sin\theta,\ -\sin\varphi\bigr).
\]
With $\mathbf{g}=\iota^{*}\overline{\mathbf{g}}$, we obtain pointwise
\[
g_{\theta\theta}
=\left\langle\partial_{\theta}\Phi,\partial_{\theta}\Phi\right\rangle_{\mathbb{R}^{3}}
=\sin^{2}\varphi,\qquad
g_{\theta\varphi}
=\left\langle\partial_{\theta}\Phi,\partial_{\varphi}\Phi\right\rangle_{\mathbb{R}^{3}}
=0,\qquad
g_{\varphi\varphi}
=\left\langle\partial_{\varphi}\Phi,\partial_{\varphi}\Phi\right\rangle_{\mathbb{R}^{3}}
=1.
\]
Therefore,
\[
(g_{ij})=
\begin{pmatrix}
\sin^{2}\varphi & 0\\
0 & 1
\end{pmatrix},
\qquad
(g^{ij})=
\begin{pmatrix}
\csc^{2}\varphi & 0\\
0 & 1
\end{pmatrix},
\qquad
\sqrt{\det(\mathbf{g})}=\sin\varphi,
\]
and the local Laplace--Beltrami formula is
\[
\Delta u
=\frac{1}{\sqrt{\det(\mathbf{g})}}  \partial_{i} \bigl(\sqrt{\det(\mathbf{g})}  g^{ij}  \partial_{j}u\bigr)=\frac{1}{\sin\varphi}  \partial_{\varphi} \bigl(\sin\varphi  \partial_{\varphi}u\bigr)
+\frac{1}{\sin^{2}\varphi}  \partial_{\theta\theta}u.
\]

Now let $N=(0,0,1)$ be the north pole and consider stereographic projection
\[
\sigma\colon S^{2}\setminus\{N\}\longrightarrow\mathbb{R}^{2},\qquad
\sigma(x,y,z)=\left(\frac{x}{1-z},\ \frac{y}{1-z}\right).
\]
Its inverse
\[
\sigma^{-1}(x,y)
=\left(\frac{2x}{1+x^{2}+y^{2}},\ \frac{2y}{1+x^{2}+y^{2}},\ \frac{x^{2}+y^{2}-1}{1+x^{2}+y^{2}}\right)
\]
composed with the inclusion $\iota\colon S^{2}\hookrightarrow\mathbb{R}^{3}$ gives
\[
\widetilde{\Phi}=\iota\circ\sigma^{-1}\colon \mathbb{R}^{2}\longrightarrow\mathbb{R}^{3},\qquad
\widetilde{\Phi}(x,y)
=\left(\frac{2x}{1+x^{2}+y^{2}},\ \frac{2y}{1+x^{2}+y^{2}},\ \frac{x^{2}+y^{2}-1}{1+x^{2}+y^{2}}\right).
\]
Differentiating and taking Euclidean inner products gives
\[
g_{xx}
=\left\langle\partial_{x}\widetilde{\Phi},\partial_{x}\widetilde{\Phi}\right\rangle_{\mathbb{R}^{3}}
=\frac{4}{\bigl(1+x^{2}+y^{2}\bigr)^{2}},\qquad
g_{yy}
=\left\langle\partial_{y}\widetilde{\Phi},\partial_{y}\widetilde{\Phi}\right\rangle_{\mathbb{R}^{3}}
=\frac{4}{\bigl(1+x^{2}+y^{2}\bigr)^{2}},
\]
\[
g_{xy}=g_{yx}
=\left\langle\partial_{x}\widetilde{\Phi},\partial_{y}\widetilde{\Phi}\right\rangle_{\mathbb{R}^{3}}
=0,
\]
whence
\[
(g_{ij})=\frac{4}{\bigl(1+x^{2}+y^{2}\bigr)^{2}}
\begin{pmatrix}
1 & 0\\
0 & 1
\end{pmatrix},
\qquad
(g^{ij})=\frac{\bigl(1+x^{2}+y^{2}\bigr)^{2}}{4}
\begin{pmatrix}
1 & 0\\
0 & 1
\end{pmatrix},
\qquad
\sqrt{\det(\mathbf{g})}=\frac{4}{\bigl(1+x^{2}+y^{2}\bigr)^{2}}.
\]
In these coordinates,
\[
\Delta u
=\frac{1}{\sqrt{\det(\mathbf{g})}}  \partial_{i} \bigl(\sqrt{\det(\mathbf{g})}  g^{ij}  \partial_{j}u\bigr)
=\frac{\bigl(1+x^{2}+y^{2}\bigr)^{2}}{4}  \bigl(\partial_{xx}u+\partial_{yy}u\bigr).
\]
\end{example}

PDEs on Riemannian manifolds arise inevitably in several fields:
\begin{itemize}
    \item In \textbf{differential geometry}, problems such as the Yamabe problem (existence of metrics of constant scalar curvature in a conformal class) and the study of geometric flows such as Ricci flow (used by Perelman in resolving the Poincaré conjecture) depend essentially on Sobolev inequalities and estimates on manifolds.
    \item In \textbf{mathematical physics}, the equations of general relativity are formulated on Lorentzian manifolds, and analysis of their associated operators requires suitable Sobolev spaces. The positive mass theorem is a paradigmatic example in which geometry, physics, and Sobolev inequalities interact.
    \item In \textbf{global analysis and topology}, the \textit{Hodge theorem} ensures that each differential cohomology class has a harmonic representative, while the \textit{Atiyah--Singer index theorem} relates the index of elliptic operators to topological invariants of the manifold. These results show how analysis of differential operators on manifolds becomes a central tool for understanding their global structure.
\end{itemize}
\begin{semblanzaHistorica}{Thierry Aubin and Sobolev spaces on manifolds}
In 1976, Thierry Aubin published a systematic treatment of Sobolev spaces on Riemannian manifolds and of the inequalities used to study geometric equations~\cite{Aubin1976Sobolev}. His work helped consolidate an intrinsic approach based on the metric and connection, particularly fruitful in conformal problems and nonlinear equations. Seeley developed spectral tools, and Strichartz extended several constructions to complete manifolds. Contemporary theory combines these ideas with uniform localization, heat kernels, and bounded geometry.
\end{semblanzaHistorica}

Sobolev spaces on manifolds are more than a generalization of Euclidean tools. They provide the language in which energy, regularity, and geometry can be studied simultaneously.

\section{Definition and basic properties}

We shall adopt an intrinsic definition using covariant derivatives. In a chart, these derivatives are expressed in terms of partial derivatives and Christoffel symbols; globally, the metric allows them to be measured and integrated without privileging coordinates.

For a function $u\in C^{\infty}(M)$, the information of order $m$ is encoded by $\nabla^{m}u$. The following results describe its local expression and its relationship with the Riemannian metric.

\begin{proposition}\label{m derivada covariante recursiva}
Let $M$ be a smooth manifold with or without boundary. Let $\nabla$ be a connection on $TM$, let $(U,\phi)$ be a chart on $M$, and let $\{\Gamma_{ij}^{k}\}$ be the Christoffel symbols of the connection with respect to $(U,\phi)$. Given $u\in C^{\infty}(M)$, $\nabla^{m}u$ has the following local expression:
\[
    (\nabla^{m}u)_{i_{1}\dots i_{m}}
    =
    \frac{\partial(\nabla^{m-1}u)_{i_{2}\dots i_{m}}}{\partial x_{i_{1}}}
    -
    \sum_{l=2}^{m}\Gamma_{i_{1}i_{l}}^{k}
    (\nabla^{m-1}u)_{i_{2}\dots i_{l-1}k i_{l+1}\dots i_{m}}.
    \]
That is, according to the convention fixed in Proposition~\ref{prop:formulas-segunda-tercera-derivada-covariante}, the new covariant index occupies the first slot.
\end{proposition}
\begin{proof}
Using Proposition~\ref{propiedades de la derivada covariante tensorial} and the convention for the position of the new index, we have
\[
\begin{aligned}
(\nabla^{m}u)_{i_{1}\dots i_{m}}
&=(\nabla_{\boldsymbol{\partial}_{i_{1}}}(\nabla^{m-1}u))
  (\boldsymbol{\partial}_{i_{2}},\dots,\boldsymbol{\partial}_{i_{m}})\\
&=\boldsymbol{\partial}_{i_{1}}\!\left(
(\nabla^{m-1}u)(\boldsymbol{\partial}_{i_{2}},\dots,
\boldsymbol{\partial}_{i_{m}})
\right)\\
&\quad-
\sum_{l=2}^{m}
(\nabla^{m-1}u)
\left(
\boldsymbol{\partial}_{i_{2}},\dots,
\boldsymbol{\partial}_{i_{l-1}},
\nabla_{\boldsymbol{\partial}_{i_{1}}}\boldsymbol{\partial}_{i_l},
\boldsymbol{\partial}_{i_{l+1}},\dots,\boldsymbol{\partial}_{i_m}
\right)\\
&=
\frac{\partial(\nabla^{m-1}u)_{i_{2}\dots i_{m}}}{\partial x_{i_{1}}}
-
\sum_{l=2}^{m}\Gamma_{i_{1}i_l}^{k}
(\nabla^{m-1}u)_{i_{2}\dots i_{l-1}k i_{l+1}\dots i_m}.
\end{aligned}
\]
\end{proof}

\begin{proposition}\label{m derivada covariante polinomio}
Let $M$ be a smooth manifold with or without boundary of dimension $n$, and let $\nabla$ be a connection on $TM$. If $u\in C^{\infty}(M)$ and $(U,\phi)$ is a smooth chart on $M$, then
$$\nabla^{m}u=\displaystyle\sum_{i_{1},\dots,i_{m}=1}^{n}\displaystyle\sum_{0\leq |\alpha|\leq m}D^{\alpha}uP^{i_{1}\dots i_{m}}_{\alpha}\mathbf{d}x^{i_{1}}\otimes \cdots\otimes \mathbf{d}x^{i_{m}}$$
on $U$, where $P^{i_{1}\dots i_{m}}_{\alpha}$ is a polynomial in the Christoffel symbols of $\nabla$ and their derivatives up to order $m-1$.
\end{proposition}
\begin{proof}

We argue by induction on $m$. For $m=1$, \(\nabla u=\displaystyle\frac{\partial u}{\partial x^{i}}\mathbf{d}x^{i}\); for $m=2$,
\[
    \nabla^{2}u
    =
    \left(
    \frac{\partial^{2}u}{\partial x^{i}\partial x^{j}}
    -\Gamma_{ij}^{k}\frac{\partial u}{\partial x^{k}}
    \right)
    \mathbf{d}x^{i}\otimes \mathbf{d}x^{j}.
    \]
Suppose the assertion holds for some $m\geq2$. By Proposition~\ref{m derivada covariante recursiva},
\[
    (\nabla^{m+1}u)_{i_{1}\dots i_{m+1}}
    =
    \frac{\partial(\nabla^{m}u)_{i_{2}\dots i_{m+1}}}{\partial x_{i_{1}}}
    -
    \sum_{l=2}^{m+1}\Gamma_{i_{1}i_{l}}^{k}
    (\nabla^{m}u)_{i_{2}\dots i_{l-1}k i_{l+1}\dots i_{m+1}}.
    \]
The induction hypothesis tells us that
\[
    (\nabla^{m}u)_{i_{2}\dots i_{m+1}}
    =
    \sum_{0\leq |\alpha|\leq m}
    D^{\alpha}uP^{i_{2}\dots i_{m+1}}_{\alpha}.
    \]
Then
\[
    \begin{aligned}
    (\nabla^{m+1}u)_{i_{1}\dots i_{m+1}}
    ={}&
    \sum_{0\leq|\alpha|\leq m}
    \frac{\partial(D^\alpha u)}{\partial x_{i_1}}
    P^{i_{2}\dots i_{m+1}}_\alpha\\
    &+
    \sum_{0\leq|\alpha|\leq m}
    D^\alpha u\,
    \frac{\partial P^{i_{2}\dots i_{m+1}}_\alpha}{\partial x_{i_1}}\\
    &-
    \sum_{l=2}^{m+1}\sum_{k=1}^{n}
    \Gamma_{i_1i_l}^{k}
    \sum_{0\leq|\alpha|\leq m}
    D^\alpha u\,
    P^{i_{2}\dots i_{l-1}k i_{l+1}\dots i_{m+1}}_\alpha .
    \end{aligned}
    \]
The first line contains derivatives of $u$ of order at most $m+1$. In the second, the coefficients are polynomials in the Christoffel symbols and their derivatives up to order $m$, because we differentiate once a polynomial that, by hypothesis, contains derivatives only up to order $m-1$. In the third line, that polynomial is multiplied by a Christoffel symbol without increasing the maximum order of derivatives of its coefficients. Grouping the terms multiplying each $D^\beta u$, $|\beta|\leq m+1$, gives the required polynomials $P^{i_1\dots i_{m+1}}_\beta$. This completes the induction.
\end{proof}

Since $(M,\mathcal{L}_{M},\lambda_{\mathbf{g}})$ is a measure space and $\lambda_{\mathbf{g}}$ is a Radon measure, the results of measure theory apply directly and allow us to consider the Lebesgue spaces $L^{p}(M)$.

If $u\in C^{\infty}(M)$, its pointwise norm is the function
\[
|\nabla^{k}u|_{\mathbf{g}}:=\sqrt{\mathbf{g}(\nabla^{k}u,\nabla^{k}u)}\colon M\longrightarrow[0,\infty).
\]
Proposition~\ref{metrica en tensores} gives its local coordinate expression:
$$|\nabla^{k}u|_{\mathbf{g}}=\sqrt{g^{i_{1}j_{1}}\cdots g^{i_{k}j_{k}}(\nabla^{k}u)_{i_{1}\dots i_{k}}(\nabla^{k}u)_{j_{1}\dots j_{k}}}.$$ Moreover, this function is at least continuous, and hence measurable. We may consider $$\int_{M}|\nabla^{k}u|_{\mathbf{g}}^{p}d\lambda_{\mathbf{g}}$$ for each $p\in[1,\infty)$. We can define $$\|u\|_{W^{m,p}(M)}:=\left(\displaystyle\sum_{j=0}^{m}\int_{M}|\nabla^{j}u|_{\mathbf{g}}^{p}d\lambda_{\mathbf{g}}\right)^{\frac{1}{p}}$$ and ask which $u\in C^{\infty}(M)$ satisfy $\|u\|_{W^{m,p}(M)}<\infty$. This motivates the following definition.
\begin{definition}\label{def: sobolev en variedades riemannianas con completaciones Hmp(M) Wmp(M) funciones escalares}
Let $(M,\mathbf{g})$ be a Riemannian manifold with or without boundary. For an integer $m\geq 0$ and $p\in [1,\infty)$, define
\begin{enumerate}[label=(\alph*)]
\item
$C^{m,p}(M):=\{u\in C^{\infty}(M)\mid\|u\|_{W^{m,p}(M)}<\infty\}$
\item $W^{m,p}(M)$ as the completion of $C^{m,p}(M)$ under the norm $\|\cdot\|_{W^{m,p}(M)}$.
\end{enumerate}
\end{definition}
\begin{remark}
We have $$C^{m,p}(M)=\{u\in C^{\infty}(M)\mid|\nabla^{j}u|_{\mathbf{g}}\in L^{p}(M)\text{ for every $j\in\{0,\dots,m\}$}\}.$$
\end{remark}
To see that this definition is consistent, we shall prove that $\left(C^{m,p}(M),\|\cdot\|_{W^{m,p}(M)}\right)$ is a normed space.
\begin{proposition}
Let $(M,\mathbf{g})$ be a Riemannian manifold with or without boundary, and let $m\in\mathbb N_0$ and $1\leq p<\infty$. Then $\left(C^{m,p}(M),\|\cdot\|_{W^{m,p}(M)}\right)$ is a normed space.
\end{proposition}
\begin{proof}
Let $u\in C^{\infty}(M)$ satisfy $\|u\|_{W^{m,p}(M)}=0$. Since $0\le\|u\|_{L^{p}(M)}\leq \|u\|_{W^{m,p}(M)}=0$, we have $\|u\|_{L^{p}(M)}=0$, so $u=0$ almost everywhere. Moreover, since $u$ is continuous, we conclude that $u\equiv 0$ on $M$.\\

Now let $\lambda\in\mathbb{R}$. Then $$\|\lambda u\|_{W^{m,p}(M)}:=\left(\displaystyle\sum_{j=0}^{m}\int_{M}|\nabla^{j}(\lambda u)|_{\mathbf{g}}^{p}d\lambda_{\mathbf{g}}\right)^{\frac{1}{p}}=\left(\displaystyle\sum_{j=0}^{m}\int_{M}|\lambda\nabla^{j}( u)|_{\mathbf{g}}^{p}d\lambda_{\mathbf{g}}\right)^{\frac{1}{p}}=\left(\displaystyle\sum_{j=0}^{m}\int_{M}|\lambda|^{p}|\nabla^{j} u|_{\mathbf{g}}^{p}d\lambda_{\mathbf{g}}\right)^{\frac{1}{p}}$$ $$=|\lambda|\left(\displaystyle\sum_{j=0}^{m}\int_{M}|\nabla^{j} u|_{\mathbf{g}}^{p}d\lambda_{\mathbf{g}}\right)^{\frac{1}{p}}=|\lambda|\|u\|_{W^{m,p}(M)},$$ where we used $|\lambda\nabla^{j}u|_{\mathbf{g}}=\sqrt{\mathbf{g}(\lambda \nabla^{j}u,\lambda\nabla^{j}u)}=\sqrt{\lambda^{2}\mathbf{g}(\nabla^{j}u,\nabla^{j}u)}=|\lambda|\sqrt{\mathbf{g}(\nabla^{j}u,\nabla^{j}u)}=|\lambda||\nabla^{j}u|_{\mathbf{g}}$ for each $j\in\{0,\dots,m\}$.\\

Let $u,v\in C^{m,p}(M)$. For $w\in C^{m,p}(M)$, define on the measure space
\[
    X_m:=M\times\{0,\dots,m\},
    \qquad
    \mu_m:=\lambda_{\mathbf{g}}\otimes\#,
    \]
the function $A_w(x,j):=|\nabla^jw(x)|_{\mathbf{g}}$, where $\#$ is counting measure. Then
\[
    \|A_w\|_{L^p(X_m,\mu_m)}
    =\left(\displaystyle\sum_{j=0}^{m}\int_M|\nabla^jw|_{\mathbf{g}}^p\,d\lambda_{\mathbf{g}}\right)^{\frac{1}{p}}
    =\|w\|_{W^{m,p}(M)}.
    \]
Since $|\nabla^j(u+v)|_{\mathbf{g}}\leq|\nabla^ju|_{\mathbf{g}}+|\nabla^jv|_{\mathbf{g}}$, we have $A_{u+v}\leq A_u+A_v$. Monotonicity of the norm and Minkowski's inequality in $L^p(X_m,\mu_m)$, given in Theorem~\ref{teo:b6-espacios-lp-espacios-de-lebesgue}, yield
\[
    \|u+v\|_{W^{m,p}(M)}
    \leq\|A_u+A_v\|_{L^p(X_m,\mu_m)}
    \leq\|u\|_{W^{m,p}(M)}+\|v\|_{W^{m,p}(M)}.
    \]

In particular, $C^{m,p}(M)$ is closed under addition and multiplication by scalars in $C^\infty(M)$; thus it is a vector space, and $\|\cdot\|_{W^{m,p}(M)}$ is a norm.
\end{proof}
We shall also define the analogue of $W^{m,p}_{0}(\Omega)$.
\begin{definition}
Let $(M,\mathbf{g})$ be a Riemannian manifold with or without boundary, and let $m\in\mathbb N_0$ and $1\leq p<\infty$. Define
\[
    W_{0}^{m,p}(M)
    :=
    \overline{C_{c}^{\infty}(\operatorname{Int}(M))}^{\,W^{m,p}(M)}.
    \]
Functions in $C_c^\infty(\operatorname{Int}(M))$ are viewed as smooth functions on $M$ by extension by zero; this is legitimate because their supports are compactly contained in the interior.
\end{definition}
The sum norm defined below is often used in the literature; see, for example, \cite{Hebey1} and \cite{Adams}. The following proposition quantifies its comparison with the norm adopted in this book.
\begin{proposition}\label{normas equivalentes de Wmp}
Let $(M,\mathbf{g})$ be a Riemannian manifold with or without boundary, and let $m\in\mathbb N_0$ and $1\leq p<\infty$. For each $u\in C^{m,p}(M)$, define
\[
\|u\|_{\widetilde{W}^{m,p}(M)}
:=\sum_{j=0}^{m}
\left(\int_M|\nabla^ju|_{\mathbf{g}}^p\,d\lambda_{\mathbf{g}}\right)^{\frac{1}{p}}.
\]
Then, for every $u\in C^{m,p}(M)$,
\[
\|u\|_{W^{m,p}(M)}
\leq\|u\|_{\widetilde{W}^{m,p}(M)}
\leq(m+1)^{1-\frac{1}{p}}\|u\|_{W^{m,p}(M)}.
\]
\end{proposition}
\begin{proof}
Set $\displaystyle x_j:=\left(\int_M|\nabla^ju|_{\mathbf{g}}^p\,d\lambda_{\mathbf{g}}\right)^{\frac{1}{p}}$.
The inequalities $\|x\|_{\ell^p(\{0,\ldots,m\})}\leq\|x\|_{\ell^1(\{0,\ldots,m\})}
\leq(m+1)^{1-\frac{1}{p}}\|x\|_{\ell^p(\{0,\ldots,m\})}$, obtained by Hölder's inequality on $\mathbb R^{m+1}$, are exactly the two asserted inequalities.
\end{proof}
Compact Riemannian manifolds play a particularly important role in the study of Sobolev spaces, since they allow generalizations of the results previously established in Euclidean spaces. If $M$ is closed, then $\operatorname{Int}(M)=M$ and $C_c^\infty(M)=C^\infty(M)$; consequently, $W_0^{m,p}(M)=W^{m,p}(M)$, and $C^\infty(M)$ is dense in this space. When $M$ has boundary, the relationship between $W_0^{m,p}(M)$ and the trace will be studied in Section~\ref{sec:sobolev-haces-frontera}, in the more general setting of sections of vector bundles.

At times, a smooth manifold $M$ may carry two Riemannian metrics, say $\mathbf{g}$ and $\widetilde{\mathbf{g}}$, making it necessary to distinguish the metric used in the Sobolev norm. In this case, we write $\|u\|_{W^{m,p}(M),\mathbf{g}}$ for the norm associated with $\mathbf{g}$ and $\|u\|_{W^{m,p}(M),\widetilde{\mathbf{g}}}$ for the norm associated with $\widetilde{\mathbf{g}}$.\\

In the compact case, other useful properties also hold; for example, the Sobolev spaces do not depend on the Riemannian metric. To prove this, we first establish the following results.
\begin{lemma}\label{lema cg}
Let $k,l\in\mathbb{N}_0$ and let $M$ be a compact smooth manifold with or without boundary equipped with two Riemannian metrics $\mathbf{g}$ and $\widetilde{\mathbf{g}}$. Then there exist $c_{k,l},C_{k,l}>0$ such that for every $p\in M$ and $v\in T^{(k,l)}(T_{p}M)$,
\[
    c_{k,l}|v|_{\mathbf g}^{2}
    \leq |v|_{\widetilde{\mathbf g}}^{2}
    \leq C_{k,l}|v|_{\mathbf g}^{2}.
    \]
Here the norms are those induced on the tensor bundle. In the proof, we shall also write $\mathbf g_p(v,v)$ and $\widetilde{\mathbf g}_p(v,v)$ for the induced inner products evaluated at $(v,v)$. For $(k,l)=(1,0)$, these are the original metrics on $T_pM$, and the conclusion is precisely the comparison in Definition~\ref{def:orden-formas-metricas}.
\end{lemma}
\begin{proof}
Let $d$ be the rank of $T^{(k,l)}(TM)$. If $M=\varnothing$ or $d=0$, the assertion is immediate. By Proposition~\ref{bolacoordenadaregular}, Remark~\ref{obs:nociones-fundamentales-bolas-coordenadas-regulares-compactas}, and compactness of $M$, there are charts $(U_i,\phi_i)$ and compact regular coordinate balls and, if $M$ has boundary, compact regular coordinate half-balls
\[
 \overline B_1,\dots,\overline B_N
\]
covering $M$ and satisfying $\overline B_i\subseteq U_i$. For each $i$, let
\[
 \Phi_i\colon T^{(k,l)}(TM)\restriction_{U_i}
 \longrightarrow U_i\times\mathbb R^d
\]
be the trivialization induced by the tensor frame of the chart containing $\overline B_i$. If $p\in\overline B_i$ and $w\in\mathbb R^d$, set
\[
 \begin{aligned}
 Q_i(p,w)
 &:=\mathbf g_p\bigl(\Phi_i^{-1}(p,w),\Phi_i^{-1}(p,w)\bigr),\\
 \widetilde Q_i(p,w)
 &:=\widetilde{\mathbf g}_p
 \bigl(\Phi_i^{-1}(p,w),\Phi_i^{-1}(p,w)\bigr).
 \end{aligned}
\]
Both functions are continuous and, for each $p$, are positive definite quadratic forms in $w$. Thus,
\[
 R_i(p,w):=\frac{\widetilde Q_i(p,w)}{Q_i(p,w)}
\]
is continuous and strictly positive on
\[
 \overline B_i\times S^{d-1},
 \qquad
 S^{d-1}:=\{w\in\mathbb R^d\mid\|w\|_2=1\}.
\]
This set is compact, so there exist
\[
 c_i:=\min_{(p,w)\in\overline B_i\times S^{d-1}}R_i(p,w)>0,
 \qquad
 C_i:=\max_{(p,w)\in\overline B_i\times S^{d-1}}R_i(p,w)<\infty.
\]
Homogeneity of the two quadratic forms implies
\[
 c_iQ_i(p,w)\leq\widetilde Q_i(p,w)\leq C_iQ_i(p,w)
\]
for every $p\in\overline B_i$ and every $w\in\mathbb R^d$. Taking
\[
 c_{k,l}:=\min_{1\leq i\leq N}c_i>0,
 \qquad
 C_{k,l}:=\max_{1\leq i\leq N}C_i<\infty,
\]
and writing $w$ for the components of $v\in T^{(k,l)}(T_pM)$ in a chart containing $p$ gives the stated inequality.
\end{proof}

The preceding argument uses only continuity and positivity of the inner products on compact subsets of trivializations. It therefore also applies on any compact set $K\subseteq M$, with continuous metrics defined on a neighborhood of $K$: simply cover $K$ by finitely many of the balls or half-balls used in the proof. The constants may depend on $K$ and on the metrics. The next proposition specifies the bounds obtained when comparison constants on the tangent bundle are already known.

\begin{proposition}[Comparison of dual norms, tensor norms, and volumes]
\label{prop:comparacion-metricas-duales-tensores-volumen}
Let $\mathbf g$ and $\widetilde{\mathbf g}$ be two Riemannian metrics on a manifold $M$ of dimension $n$, and let $A\subseteq M$. Suppose that $0<c\leq C<\infty$ and
\[
 c\,\mathbf g_p(v,v)\leq\widetilde{\mathbf g}_p(v,v)
 \leq C\,\mathbf g_p(v,v)
 \qquad(p\in A,\ v\in T_pM).
\]
Then the dual metrics satisfy
\begin{equation}
\label{eq:comparacion-duales-cuantitativa}
 C^{-1}\mathbf g_p^{-1}(\xi,\xi)
 \leq\widetilde{\mathbf g}_p^{-1}(\xi,\xi)
 \leq c^{-1}\mathbf g_p^{-1}(\xi,\xi)
 \qquad(\xi\in T_p^*M,\ p\in A).
\end{equation}
Moreover, for $k,l\in\mathbb N_0$ and $\mathbf T\in T^{(k,l)}(T_pM)$,
\begin{equation}
\label{eq:comparacion-tensores-cuantitativa}
 c^{k/2}C^{-l/2}|\mathbf T|_{\mathbf g}
 \leq|\mathbf T|_{\widetilde{\mathbf g}}
 \leq C^{k/2}c^{-l/2}|\mathbf T|_{\mathbf g}.
\end{equation}
In particular, a comparison $K^{-1}\mathbf g\leq\widetilde{\mathbf g}\leq K\mathbf g$, with $K\geq1$, gives the constants $K^{-(k+l)/2}$ and $K^{(k+l)/2}$ for the norms of tensors of type $(k,l)$.

The ratio of the volume elements is the positive function
\[
 \rho(p)=
 \frac{\sqrt{\det(\widetilde{\mathbf g}_p)}}
      {\sqrt{\det(\mathbf g_p)}},
\]
where both determinants are computed in the same chart. This function is independent of the chart, satisfies $d\lambda_{\widetilde{\mathbf g}}=\rho\,d\lambda_{\mathbf g}$, and obeys
\begin{equation}
\label{eq:comparacion-volumen-cuantitativa}
 c^{n/2}\leq\rho(p)\leq C^{n/2}
 \qquad(p\in A).
\end{equation}
Consequently, if $A$ is measurable and $f\colon A\to[0,\infty]$ is measurable, then
\[
 c^{n/2}\int_A f\,d\lambda_{\mathbf g}
 \leq\int_A f\,d\lambda_{\widetilde{\mathbf g}}
 \leq C^{n/2}\int_A f\,d\lambda_{\mathbf g}.
\]
\end{proposition}

\begin{proof}
Fix $p\in A$. There is a unique endomorphism $L_p\colon T_pM\to T_pM$ such that
\[
 \widetilde{\mathbf g}_p(v,w)=\mathbf g_p(L_pv,w)
 \qquad(v,w\in T_pM).
\]
It is obtained by applying the inverse of the isomorphism $w\mapsto\mathbf g_p(w,\cdot)$ to the covector $\widetilde{\mathbf g}_p(v,\cdot)$. Symmetry of $\widetilde{\mathbf g}_p$ implies that $L_p$ is self-adjoint with respect to $\mathbf g_p$, and its positivity implies $\mathbf g_p(L_pv,v)>0$ when $v\neq0$. The spectral theorem provides a basis $e_1,\ldots,e_n$ orthonormal for $\mathbf g_p$ and numbers $\lambda_1,\ldots,\lambda_n>0$ such that
\[
 L_pe_i=\lambda_i e_i,\qquad
 \mathbf g_p(e_i,e_j)=\delta_{ij},\qquad
 \widetilde{\mathbf g}_p(e_i,e_j)=\lambda_i\delta_{ij}.
\]
Applying the hypothesis to $v=e_i$ gives $c\leq\lambda_i\leq C$ for every $i$. This basis is chosen in the fixed fiber $T_pM$; the chosen bases need not vary smoothly with $p$.

Let $e^1,\ldots,e^n$ be the dual basis. The matrices of the dual metrics are the inverses of the matrices of the tangent metrics. Thus, for $\xi=\displaystyle\sum_{i=1}^{n}\xi_i e^i$,
\[
 \mathbf g_p^{-1}(\xi,\xi)=\sum_{i=1}^{n}\xi_i^2,
 \qquad
 \widetilde{\mathbf g}_p^{-1}(\xi,\xi)
 =\sum_{i=1}^{n}\lambda_i^{-1}\xi_i^2.
\]
Multiplying the inequalities $C^{-1}\leq\lambda_i^{-1}\leq c^{-1}$ by $\xi_i^2$ and summing proves \eqref{eq:comparacion-duales-cuantitativa}. This explains why the lower and upper constants exchange roles when comparing dual metrics.

Now consider the basis tensors
\[
 E_{I,J}:=
 e_{i_1}\otimes\cdots\otimes e_{i_k}
 \otimes e^{j_1}\otimes\cdots\otimes e^{j_l},
\]
where $I=(i_1,\ldots,i_k)$ and $J=(j_1,\ldots,j_l)$ range over the indices from $1$ to $n$. They are orthogonal for both tensor inner products, and
\[
 |E_{I,J}|_{\mathbf g}^{2}=1,\qquad
 |E_{I,J}|_{\widetilde{\mathbf g}}^{2}
 =\lambda_{i_1}\cdots\lambda_{i_k}
  \lambda_{j_1}^{-1}\cdots\lambda_{j_l}^{-1}.
\]
An empty product is interpreted as $1$. Each weight in the last expression lies between $c^kC^{-l}$ and $C^kc^{-l}$. If $\mathbf T=\displaystyle\sum_{\substack{I\in\{1,\dots,n\}^k\\J\in\{1,\dots,n\}^l}}T^{I}_{J}E_{I,J}$, orthogonality gives
\[
 |\mathbf T|_{\mathbf g}^{2}=\sum_{\substack{I\in\{1,\dots,n\}^k\\J\in\{1,\dots,n\}^l}}|T^{I}_{J}|^2,\qquad
 |\mathbf T|_{\widetilde{\mathbf g}}^{2}
 =\sum_{\substack{I\in\{1,\dots,n\}^k\\J\in\{1,\dots,n\}^l}}|T^{I}_{J}|^2
   \lambda_{i_1}\cdots\lambda_{i_k}
   \lambda_{j_1}^{-1}\cdots\lambda_{j_l}^{-1}.
\]
Multiplying the bounds for each weight by $|T^{I}_{J}|^2$ and summing yields
\[
 c^kC^{-l}|\mathbf T|_{\mathbf g}^{2}
 \leq|\mathbf T|_{\widetilde{\mathbf g}}^{2}
 \leq C^kc^{-l}|\mathbf T|_{\mathbf g}^{2}.
\]
Taking square roots proves \eqref{eq:comparacion-tensores-cuantitativa}.

In the basis $e_1,\ldots,e_n$, the ratio of determinants is $\lambda_1\cdots\lambda_n$. A simultaneous change of basis multiplies both determinants by the same square of the determinant of the change-of-basis matrix, leaving their ratio unchanged. In particular,
\[
 \frac{\det(\widetilde{\mathbf g}_p)}{\det(\mathbf g_p)}
 =\prod_{i=1}^n\lambda_i,\qquad
 c^n\leq\prod_{i=1}^n\lambda_i\leq C^n.
\]
The coordinate formula for Riemannian measure implies $d\lambda_{\widetilde{\mathbf g}}=\rho\,d\lambda_{\mathbf g}$, and the square roots of the last bounds prove \eqref{eq:comparacion-volumen-cuantitativa}. Finally,
\[
 \int_A f\,d\lambda_{\widetilde{\mathbf g}}
 =\int_A f\rho\,d\lambda_{\mathbf g}.
\]
Multiply \eqref{eq:comparacion-volumen-cuantitativa} by $f\geq0$ and use monotonicity of the integral. The argument also holds when one of the integrals is infinite.
\end{proof}

This proposition compares two norms of the same tensor. In studying Sobolev norms, we must also compare the tensors $\nabla^j u$ and $\widetilde\nabla^j u$ when the connection changes. For functions, $\nabla u=\widetilde\nabla u=du$; from order two onward, covariant derivatives may differ. The following lemma compares the full collections of derivatives up to order $m$, including their lower-order terms.
\begin{lemma}\label{lema dificil de probar}
Let $m\geq 0$ and let $M$ be a compact smooth manifold with or without boundary, with two Riemannian metrics $\mathbf{g}$ and $\widetilde{\mathbf{g}}$, and let $\nabla$ and $\widetilde{\nabla}$ be two connections on $TM$. Then there exist $c_{m},C_{m}>0$ such that $$c_{m}\left(\displaystyle\sum_{j=0}^{m}|\nabla^{j}u|^{2}_{\mathbf{g}}\right)\leq \displaystyle\sum_{j=0}^{m}|\widetilde{\nabla}^{j}u|^{2}_{\widetilde{\mathbf{g}}}\leq C_{m}\left(\displaystyle\sum_{j=0}^{m}|\nabla^{j}u|^{2}_{\mathbf{g}}\right)$$ for every $u\in C^{\infty}(M)$.
\end{lemma}
\begin{proof}
If $M=\varnothing$, the assertion is immediate. Assume $M\neq\varnothing$. By Definition~\ref{suma de whitney}, we may consider the bundle
\[
 \mathscr T_m:=\bigoplus_{j=0}^{m}T^{(0,j)}(TM),
 \qquad T^{(0,0)}(TM):=M\times\mathbb R.
\]
For $u\in C^\infty(M)$, write
\[
 J_\nabla^m u:=(u,\nabla u,\dots,\nabla^m u)
\]
and, for each $p\in M$, define
\[
 \mathscr D_{\nabla,p}
 :=\{J_\nabla^m u(p)\mid u\in C^\infty(M)\}
 \subseteq(\mathscr T_m)_p.
\]
Linearity of $\nabla^j$ implies that $\mathscr D_{\nabla,p}$ is a vector subspace of $(\mathscr T_m)_p$.

Let $(U,\phi)$ be a smooth chart, let $p\in U$, and set $x:=\phi(p)$ and
\[
 N:=\sum_{j=0}^{m}{n+j-1\choose j}={n+m\choose m}.
\]
With the multi-index notation from Definition~\ref{def:nociones-fundamentales-multiindice}, define the linear map
\[
 F_{\nabla,p}^{\phi}\colon\mathscr D_{\nabla,p}\longrightarrow\mathbb R^N,
 \qquad
 F_{\nabla,p}^{\phi}(J_\nabla^m u(p))
 :=\bigl(D^\alpha(u\circ\phi^{-1})(x)\bigr)_{|\alpha|\leq m}.
\]
Propositions~\ref{m derivada covariante recursiva} and \ref{m derivada covariante polinomio} give, for $j\in\{1,\dots,m\}$,
\begin{equation}
 (\nabla^ju)_{i_1\dots i_j}
 =\partial_{i_1}\cdots\partial_{i_j}(u\circ\phi^{-1})
 +\sum_{|\alpha|<j}
 A^\alpha_{i_1\dots i_j}
 D^\alpha(u\circ\phi^{-1}),
 \label{eq:jet-covariante-triangular-funciones}
\end{equation}
where the coefficients $A^\alpha_{i_1\dots i_j}$ are smooth on $U$. By induction on $j$, the equality $J_\nabla^m u(p)=J_\nabla^m v(p)$ determines first $u(p)=v(p)$ and then all partial derivatives of $u\circ\phi^{-1}$ and $v\circ\phi^{-1}$ up to order $m$. Thus $F_{\nabla,p}^{\phi}$ is well defined and injective.

Let $a=(a_\alpha)_{|\alpha|\leq m}\in\mathbb R^N$. The polynomial
\[
 P_a(y):=\sum_{|\alpha|\leq m}
 \frac{a_\alpha}{\alpha!}(y-x)^\alpha
\]
satisfies $D^\alpha P_a(x)=a_\alpha$ for every $|\alpha|\leq m$. Theorem~\ref{extension de funciones suaves}, applied to $P_a\circ\phi$, provides a function $u\in C^\infty(M)$ agreeing with $P_a\circ\phi$ on a neighborhood of $p$. Then $F_{\nabla,p}^{\phi}(J_\nabla^m u(p))=a$, and $F_{\nabla,p}^{\phi}$ is an isomorphism.

Let $(e_1,\dots,e_N)$ be the standard basis of $\mathbb R^N$. Substituting each $e_\nu$ as the prescribed partial jet in \eqref{eq:jet-covariante-triangular-funciones} shows that the components of
\[
 \boldsymbol\sigma_\nu(q)
 :=(F_{\nabla,q}^{\phi})^{-1}(e_\nu),
 \qquad q\in U,\quad \nu\in\{1,\dots,N\},
\]
in the coordinate tensor frames are smooth functions. Moreover, $\boldsymbol\sigma_1(q),\dots,\boldsymbol\sigma_N(q)$ form a basis of $\mathscr D_{\nabla,q}$ for each $q\in U$. Lemma~\ref{criterio del marco local suave para subhaces} implies that
\[
 \mathscr D_\nabla:=\bigcup_{p\in M}\mathscr D_{\nabla,p}
\]
is a smooth vector subbundle of $\mathscr T_m$ of rank $N$. The same construction applied to $\widetilde\nabla$ produces the subbundle $\mathscr D_{\widetilde\nabla}$ and the isomorphisms $F_{\widetilde\nabla,p}^{\phi}\colon
\mathscr D_{\widetilde\nabla,p}\longrightarrow\mathbb R^N$.

Let $(\overline B,\phi)$ be a compact regular coordinate ball or half-ball whose coordinate map is the restriction of the preceding chart. For $p\in\overline B$ and $v\in\mathbb R^N$, define
\[
 \begin{aligned}
 Q(p,v)
 &:=\displaystyle\sum_{j=0}^{m}
 \left|\bigl((F_{\nabla,p}^{\phi})^{-1}(v)\bigr)_j
 \right|_{\mathbf g}^{2},\\
 \widetilde Q(p,v)
 &:=\displaystyle\sum_{j=0}^{m}
 \left|\bigl((F_{\widetilde\nabla,p}^{\phi})^{-1}(v)\bigr)_j
 \right|_{\widetilde{\mathbf g}}^{2}.
 \end{aligned}
\]
The families $\boldsymbol\sigma_1,\dots,\boldsymbol\sigma_N$ constructed above and their analogues for $\widetilde\nabla$ show that $Q$ and $\widetilde Q$ are continuous on $\overline B\times\mathbb R^N$. For each $p$, both are positive definite quadratic forms in $v$. Thus the ratio
\[
 R(p,v):=\frac{\widetilde Q(p,v)}{Q(p,v)}
\]
is continuous and strictly positive on $\overline B\times S^{N-1}$, where $S^{N-1}:=\{v\in\mathbb R^N\mid\|v\|_2=1\}$. Compactness of $\overline B\times S^{N-1}$ provides numbers
\[
 c_{\overline B,m}
 :=\min_{(p,v)\in\overline B\times S^{N-1}}R(p,v)>0,
 \qquad
 C_{\overline B,m}
 :=\max_{(p,v)\in\overline B\times S^{N-1}}R(p,v)<\infty.
\]
Since $Q$ and $\widetilde Q$ are homogeneous of degree two in $v$, it follows that
\begin{equation}
 c_{\overline B,m}Q(p,v)
 \leq\widetilde Q(p,v)
 \leq C_{\overline B,m}Q(p,v)
 \label{eq:comparacion-local-formas-cuadraticas-jets}
\end{equation}
for every $p\in\overline B$ and every $v\in\mathbb R^N$; for $v\neq0$, apply the definitions of $c_{\overline B,m}$ and $C_{\overline B,m}$ to $v/\|v\|_2$, and for $v=0$ the inequality is immediate.

For $p\in\overline B$ and $u\in C^\infty(M)$, the definition of the two isomorphisms gives
\[
 F_{\nabla,p}^{\phi}(J_\nabla^m u(p))
 =\bigl(D^\alpha(u\circ\phi^{-1})(\phi(p))\bigr)_{|\alpha|\leq m}
 =F_{\widetilde\nabla,p}^{\phi}(J_{\widetilde\nabla}^m u(p)).
\]
Substituting this vector into \eqref{eq:comparacion-local-formas-cuadraticas-jets} yields
\[
 c_{\overline B,m}\sum_{j=0}^{m}|\nabla^ju(p)|_{\mathbf g}^{2}
 \leq
 \sum_{j=0}^{m}|\widetilde\nabla^ju(p)|_{\widetilde{\mathbf g}}^{2}
 \leq
 C_{\overline B,m}\sum_{j=0}^{m}|\nabla^ju(p)|_{\mathbf g}^{2}.
\]

Compactness of $M$ allows us to choose a finite cover $\overline B_1,\dots,\overline B_k$ by compact regular coordinate balls or half-balls. Taking
\[
 c_m:=\min_{1\leq i\leq k}c_{\overline B_i,m}>0,
 \qquad
 C_m:=\max_{1\leq i\leq k}C_{\overline B_i,m}<\infty,
\]
gives the two-sided inequality on all of $M$.
\end{proof}
We would also like a way to relate integrals with respect to two metrics when $M$ is compact; this can be obtained using Lemma \ref{lema cg}.
\begin{lemma}\label{desigualdad integrales}
Let $M$ be a compact smooth manifold with or without boundary equipped with two Riemannian metrics $\mathbf{g}$ and $\widetilde{\mathbf{g}}$. Then there exist $C,c>0$ such that for every measurable $f\colon M\longrightarrow [0,\infty]$, $$c\int_{M}fd\lambda_{\mathbf{g}}\leq\int_{M}fd\lambda_{\widetilde{\mathbf{g}}}\leq C\int_{M}fd\lambda_{\mathbf{g}}. $$
\end{lemma}
\begin{proof}
Consider a finite family $(\overline{B}_{i},\phi_{i})_{i=1}^{k}$ of compact regular coordinate balls and, where appropriate, half-balls whose relative interiors cover $M$ (this is possible because $M$ is compact). We have $\sqrt{\det(\mathbf{g})},\sqrt{\det(\widetilde{\mathbf{g}})}>0$ on $\overline{B}_{i}$ for every $i\in \{1,\dots,k\}$, since the matrices of both metrics are symmetric and positive definite. We may therefore define $h_{i}\colon \overline{B}_{i}\longrightarrow [0,\infty)$, $h_{i}(p)=\frac{\sqrt{\det(\widetilde{\mathbf{g}}_{p})}}{\sqrt{\det(\mathbf{g}_{p})}}$. Note that $h_{i}>0$. Since $\overline{B}_{i}$ is compact, there exist $c_{i}=\displaystyle\min_{p\in \overline{B}_{i}}h_i(p)$ and $C_{i}=\displaystyle\max_{p\in \overline{B}_{i}}h_i(p)$, both positive. Then $c_{i}\sqrt{\det(\mathbf{g})}\leq\sqrt{\det(\widetilde{\mathbf{g}})}\leq C_{i}\sqrt{\det(\mathbf{g})} $ on $\overline{B}_{i}$. Let $(A_{i})_{i=1}^{k}$ be a family of pairwise disjoint measurable sets as in Definition \ref{medibles ajenos dos a dos, definicion de medida Riemann--Lebesgue} for $A=M$. We have $M=\displaystyle\coprod_{i=1}^{k}A_{i}$. Let $c:=\displaystyle\min_{1\leq i\leq k}c_{i}$ and $C:=\displaystyle\max_{1\leq i\leq k}C_{i}$. We have $f=\displaystyle\sum_{i=1}^{k}f\mathbf{1}_{A_{i}}$. Moreover, $$c_{i}\displaystyle\int
_{A_{i}}fd\lambda_{\mathbf{g}}=c_i\int_{\phi_{i}(A_{i})}(f\circ \phi_{i}^{-1})\sqrt{\det(\mathbf{g})\circ \phi_{i}^{-1}}d\lambda_{n}$$ $$\leq \int_{\phi_{i}(A_{i})}(f\circ \phi_{i}^{-1})\sqrt{\det(\widetilde{\mathbf{g}})\circ \phi_{i}^{-1}}d\lambda_{n}=\int_{A_{i}}fd\lambda_{\widetilde{\mathbf{g}}}.$$ Similarly, $$\int_{A_{i}}fd\lambda_{\widetilde{\mathbf{g}}}\leq C_{i} \int_{A_{i}}fd\lambda_{\mathbf{g}}.$$ Therefore, $$c_{i}\displaystyle\int_{M}f\mathbf{1}_{A_{i}}d\lambda_{\mathbf{g}}\leq \int_{M}f\mathbf{1}_{A_{i}}d\lambda_{\widetilde{\mathbf{g}}}\leq C_{i} \int_{M}f\mathbf{1}_{A_{i}}d\lambda_{\mathbf{g}}.$$ Hence $$c\int_{M}fd\lambda_{\mathbf{g}}=c\displaystyle\sum_{i=1}^{k}\int_{M}f\mathbf{1}_{A_{i}}d\lambda_{\mathbf{g}}\leq \displaystyle\sum_{i=1}^{k}c_{i}\displaystyle\int_{M}f\mathbf{1}_{A_{i}}d\lambda_{\mathbf{g}}\leq \displaystyle\sum_{i=1}^{k}\int_{M}f\mathbf{1}_{A_{i}}d\lambda_{\widetilde{\mathbf{g}}}=\int_{M}fd\lambda_{\widetilde{\mathbf{g}}}.$$
Similarly, $$\int_{M}fd\lambda_{\widetilde{\mathbf{g}}}\leq C\int_{M}fd\lambda_{\mathbf{g}}.$$

\end{proof}
\begin{remark}\label{desigualdad integrales triple}
If $u,v,w\colon M\longrightarrow [0,\infty]$ are measurable and satisfy $u\leq v\leq w$, then Lemma \ref{desigualdad integrales} gives $$c\int_{M}ud\lambda_{\mathbf{g}}\leq \int_{M}ud\lambda_{\widetilde{\mathbf{g}}}\leq \int_{M}vd\lambda_{\widetilde{\mathbf{g}}}\leq \int_{M}wd\lambda_{\widetilde{\mathbf{g}}}\leq C\int_{M}wd\lambda_{\mathbf{g}}, $$; that is, $$c\int_{M}ud\lambda_{\mathbf{g}}\leq \int_{M}vd\lambda_{\widetilde{\mathbf{g}}}\leq C\int_{M}wd\lambda_{\mathbf{g}}.$$
\end{remark}
\begin{theorem}\label{W no depende de metrica}
Let \(M\) be a compact smooth manifold with or without boundary, let \(m\in\mathbb N_0\), and let \(1\leq p<\infty\). If \(\mathbf g\) and \(\widetilde{\mathbf g}\) are two Riemannian metrics on \(M\), then the norms \(\|\cdot\|_{W^{m,p}(M),\mathbf g}\) and \(\|\cdot\|_{W^{m,p}(M),\widetilde{\mathbf g}}\) are equivalent. In particular, the two completions canonically define the same space \(W^{m,p}(M)\).

\end{theorem}
\begin{proof}
Let \(\nabla\) and \(\widetilde\nabla\) be the Levi--Civita connections of \(\mathbf g\) and \(\widetilde{\mathbf g}\), respectively. For \(u\in C^\infty(M)\) and \(q\in M\), set
\[
x(q):=\bigl(|u(q)|,|\nabla u(q)|_{\mathbf g},\dots,
|\nabla^m u(q)|_{\mathbf g}\bigr),
\]
\[
y(q):=\bigl(|u(q)|,|\widetilde\nabla u(q)|_{\widetilde{\mathbf g}},
\dots,|\widetilde\nabla^m u(q)|_{\widetilde{\mathbf g}}\bigr).
\]
By Lemma~\ref{lema dificil de probar}, there exist \(c_m,C_m>0\) such that
\[
\sqrt{c_m}\,\|x(q)\|_2
\leq\|y(q)\|_2
\leq\sqrt{C_m}\,\|x(q)\|_2
\qquad(q\in M).
\]
Since all norms on \(\mathbb R^{m+1}\) are equivalent, there exist \(a_{m,p},b_{m,p}>0\) such that
\[
a_{m,p}\|z\|_p\leq\|z\|_2\leq b_{m,p}\|z\|_p
\qquad(z\in\mathbb R^{m+1}).
\]
Applying these two inequalities first to \(x(q)\) and then to \(y(q)\) yields, pointwise,
\[
\frac{a_{m,p}\sqrt{c_m}}{b_{m,p}}\|x(q)\|_p
\leq\|y(q)\|_p
\leq
\frac{b_{m,p}\sqrt{C_m}}{a_{m,p}}\|x(q)\|_p.
\]
Raise this to the power \(p\). Lemma~\ref{desigualdad integrales}, applied to the two Riemannian measures, allows us to integrate the outer expressions and gives constants \(A,B>0\), independent of \(u\), for which
\[
A\|u\|_{W^{m,p}(M),\mathbf g}
\leq
\|u\|_{W^{m,p}(M),\widetilde{\mathbf g}}
\leq
B\|u\|_{W^{m,p}(M),\mathbf g}.
\]
Thus, the identity on \(C^\infty(M)\) is bicontinuous for the two norms. By the universal property of completion, it extends uniquely to a bicontinuous isomorphism between the two completions, identifying each smooth function with itself.
\end{proof}

\begin{example}[Dependence of the constants on derivatives of the metric]
\label{ej:comparacion-cero-no-controla-sobolev-superior}
Comparison of metrics at order zero does not suffice to control the equivalence constants of all Sobolev norms. Consider the compact manifold with boundary $M=[0,2\pi]$, with $\mathbf g_0=(\mathbf dx)^2$, and, for each integer $N\geq1$, the metric
\[
 \mathbf g_N=e^{2\sin(Nx)}(\mathbf dx)^2.
\]
For every $x\in M$ and every $v\in T_xM$,
\[
 e^{-2}\mathbf g_0(v,v)\leq\mathbf g_N(v,v)
 \leq e^2\mathbf g_0(v,v).
\]
The constants are independent of $N$. However, the only coefficient of the Levi--Civita connection of $\mathbf g_N$ is
\[
 \Gamma^1_{11}(\mathbf g_N)
 =\frac12 e^{-2\sin(Nx)}
       \frac{d}{dx}\bigl(e^{2\sin(Nx)}\bigr)
 =N\cos(Nx).
\]
For the fixed smooth function $u(x)=x$, we have $u'=1$, $u''=0$, and, by the coordinate formula for the Hessian,
\[
 (\nabla^{\mathbf g_N})^2u
 =\bigl(u''-\Gamma^1_{11}(\mathbf g_N)u'\bigr)(\mathbf dx)^2
 =-N\cos(Nx)(\mathbf dx)^2.
\]
Since $|(\mathbf dx)^2|_{\mathbf g_N}^{2}=e^{-4\sin(Nx)}$ and $d\lambda_{\mathbf g_N}=e^{\sin(Nx)}\,dx$, we obtain
\[
 \begin{aligned}
 \|(\nabla^{\mathbf g_N})^2u\|_{L^2(M,\mathbf g_N)}^2
 &=N^2\int_0^{2\pi}\cos^2(Nx)e^{-3\sin(Nx)}\,dx\\
 &\geq e^{-3}N^2\int_0^{2\pi}\cos^2(Nx)\,dx
 =\pi e^{-3}N^2.
 \end{aligned}
\]
Therefore, $\|u\|_{W^{2,2}(M),\mathbf g_N}$ tends to infinity, while $\|u\|_{W^{2,2}(M),\mathbf g_0}$ remains fixed. Theorem~\ref{W no depende de metrica} does give equivalence for each fixed pair $(\mathbf g_0,\mathbf g_N)$. What this example rules out is a common constant deduced solely from the bounds $e^{-2}$ and $e^2$. The additional control of the connections and their derivatives appears in Lemma~\ref{lema dificil de probar}.
\end{example}

\begin{remark}[Comparison on compact coordinate subsets]
\label{obs:comparacion-metricas-compactos-coordenados}
The comparisons used in the preceding proof are pointwise on the compact coordinate subsets: they are obtained by bounding the ratio of two positive quadratic forms on a compact sphere and the ratio of the volume elements. They therefore also hold on a coordinate ball or half-ball $U$ whose closure is contained in an ambient chart, for smooth metrics and connections defined on that chart. It suffices to apply the bounds on $\overline U$ and integrate over $U$; $\overline U$ need not be a smooth manifold with boundary. This matters for half-balls, whose closures have a corner where the flat and spherical portions meet. In the following formulas, norms of smooth functions on $\overline U$ denote these same integrals: the coordinate boundary has zero $n$-dimensional measure.
\end{remark}

\subsection{Concrete Sobolev functions on a manifold}

The scalar object is still a function on $M$, as in Chapter~1.
We now consider classes of functions equal almost everywhere with respect
to $\lambda_{\mathbf g}$ and require integrability of their derivatives.
The metric measures the covectors $du$ and specifies the volume element;
a chart permits calculation. The following examples connect the
completion definition with functions that can be described point by point.

\begin{example}[Height and a corner along the equator]
\label{ej:sobolev-altura-y-esquina-esfera}
On the unit sphere with its induced metric, use $h(p)=z$ from
Example~\ref{ej:altura-esfera-cartas-autocontenido}.
Since $S^2$ is compact, every covariant derivative of the smooth function
$h$ is bounded and the volume is finite. Hence
$h\in W^{m,p}(S^2)$ for all $m\in\mathbb N_0$ and $1\leq p<\infty$.

Its first-order norm can be computed explicitly. For $w\in T_pS^2$,
one has $dh_p(w)=w_3$ and
\[
 \operatorname{grad}h(p)=e_3-zp,
 \qquad |\operatorname{grad}h(p)|^2=1-z^2,
\]
since $e_3-zp$ is tangent and its inner product with $w$ is $w_3$.
In coordinates
\[
 F(\theta,z)
 =(\sqrt{1-z^2}\cos\theta,\sqrt{1-z^2}\sin\theta,z),
 \quad 0<\theta<2\pi,\quad -1<z<1,
\]
the vectors $F_\theta,F_z$ are orthogonal and
\[
 |F_\theta|^2=1-z^2,\qquad |F_z|^2=(1-z^2)^{-1}.
\]
Thus the metric determinant is $1$, and
$d\lambda_{\mathbf g}=d\theta\,dz$ in this chart.
The omitted meridian and poles have measure zero, as follows in charts
from the null measure of lines and points in $\mathbb R^2$.
Therefore,
\[
 \int_{S^2}h^2\,d\lambda_{\mathbf g}
 =2\pi\int_{-1}^1z^2\,dz=\frac{4\pi}{3},
 \qquad
 \int_{S^2}|dh|_{\mathbf g}^2\,d\lambda_{\mathbf g}
 =2\pi\int_{-1}^1(1-z^2)\,dz=\frac{8\pi}{3}.
\]
In particular, $\|h\|_{W^{1,2}(S^2)}^2=4\pi$.

The same scalar datum provides a PDE example on the sphere.
The coordinate Laplacian formula, with the determinant just computed,
gives
\[
 \Delta_{\mathbf g}h
 =\partial_\theta\bigl((1-z^2)^{-1}\partial_\theta z\bigr)
  +\partial_z\bigl((1-z^2)\partial_z z\bigr)
 =-2z=-2h.
\]
Both sides are smooth, so continuity extends the identity globally,
including the points omitted by the chart.
Integration by parts on the boundaryless sphere yields
\[
 \int_{S^2}\langle dh,d\varphi\rangle_{\mathbf g}\,d\lambda_{\mathbf g}
 =2\int_{S^2}h\varphi\,d\lambda_{\mathbf g},
 \qquad \varphi\in C^\infty(S^2).
\]
Thus $-\Delta_{\mathbf g}h=2h$ compares scalar functions on the same
domain, and its weak formulation measures and integrates their
differentials defined on tangent spaces without requiring global coordinates.

The function $u:=|h|$ is continuous and is not $C^1$ on the equator,
as proved in Example~\ref{ej:regularidad-escalar-esfera}. Nevertheless,
$u\in W^{1,p}(S^2)$ for every $1\leq p<\infty$.
To justify this using the completion definition, consider
\[
 u_\varepsilon:=\sqrt{h^2+\varepsilon^2},\qquad
 du_\varepsilon=\frac{h}{\sqrt{h^2+\varepsilon^2}}\,dh,
 \qquad 0<\varepsilon<1.
\]
The functions $u_\varepsilon$ are smooth and
$0\leq u_\varepsilon-u\leq\varepsilon$.
Set $v:=\operatorname{sgn}(h)\,dh$, assigning zero on the equator.
Since $|dh|_{\mathbf g}\leq1$, for $0<a<1$ the difference
$|du_\varepsilon-v|_{\mathbf g}$ is at most $1$ where $|h|<a$.
Where $|h|\geq a$,
\[
 0\leq1-\frac{|h|}{\sqrt{h^2+\varepsilon^2}}
 =\frac{\varepsilon^2}
 {\sqrt{h^2+\varepsilon^2}(\sqrt{h^2+\varepsilon^2}+|h|)}
 \leq\frac{\varepsilon^2}{2a^2}.
\]
The band $|h|<a$ has area $4\pi a$ by the preceding coordinates.
Taking $a=\sqrt{\varepsilon}$ yields the bounds
\[
 \|u_\varepsilon-u\|_{L^p(S^2)}^p\leq4\pi\varepsilon^p,
 \qquad
 \|du_\varepsilon-v\|_{L^p(S^2)}^p
 \leq4\pi\sqrt{\varepsilon}+4\pi(\varepsilon/2)^p.
\]
Thus $(u_\varepsilon)$ is Cauchy in the $W^{1,p}$ norm, its $L^p$ limit
is $u$, and the limit of its differentials is $v$.
The completion therefore contains the class of $u$, with weak
differential $du=\operatorname{sgn}(h)\,dh$.
The weak identity also follows by passing to the limit in integration
by parts for $u_\varepsilon$.
Since $u^2=h^2$ and $|du|_{\mathbf g}^2=|dh|_{\mathbf g}^2$ away from
the equator, $\|u\|_{W^{1,2}(S^2)}^2=4\pi$.
The corner does not destroy integrability of the differential.

Moreover, setting $\varepsilon=1/k$, this sequence does not converge
in the Sobolev norm to a smooth function. If $w$ were such a limit,
it would equal $|h|$ almost everywhere by convergence in $L^p$.
Continuity of both functions and positivity of the measure of every
nonempty open set imply $w=|h|$ at every point.
But for $\gamma(t)=(\sqrt{1-t^2},0,t)$, $|t|<1/2$, we have
$w(\gamma(t))=|t|$, which is not differentiable at zero.
Composition of a smooth function with a smooth curve would be
smooth, a contradiction. In particular, $C^{1,p}(S^2)$ is
incomplete in this norm.
\end{example}

\begin{example}[Incompleteness of smooth functions on the sphere at every order]
\label{ej:incompletitud-sobolev-esfera-orden-arbitrario}
Let $S^2$ be the unit sphere with its induced metric, and let
$h(p)=z$ be its height function. Fix $m\in\mathbb N_0$ and
$1\leq p<\infty$. We construct a sequence in $C^{m,p}(S^2)$
that is Cauchy in the Sobolev norm and does not converge in
that norm to any smooth function.

Define on $\mathbb R$
\[
 H_0(t):=\mathbf1_{(0,\infty)}(t),\qquad
 H_r(t):=\frac{(t_+)^r}{r!}\quad(r\geq1),\qquad
 t_+:=\max\{t,0\}.
\]
For $r\geq1$ and $\varphi\in C_c^\infty(\mathbb R)$,
\[
 \begin{aligned}
 \int_{\mathbb R}H_r(t)\varphi'(t)\,dt
 &=\int_0^\infty\frac{t^r}{r!}\varphi'(t)\,dt\\
 &=\left[\frac{t^r}{r!}\varphi(t)\right]_0^\infty
          -\int_0^\infty\frac{t^{r-1}}{(r-1)!}\varphi(t)\,dt\\
 &=-\int_{\mathbb R}H_{r-1}(t)\varphi(t)\,dt.
 \end{aligned}
\]
The term at zero vanishes because $r\geq1$, and that at infinity
because the test function has compact support. Iteration gives
$D^jH_m=H_{m-j}$ for $0\leq j\leq m$ as distributions;
for $m=0$ only the order-zero equality is used.

Take a nonnegative mollifier
$\rho\in C_c^\infty((-1,1))$ of integral one, and set
\[
 \varepsilon_k:=\frac{1}{k+1},\qquad
 a_k:=\rho_{\varepsilon_k}*H_m,
 \qquad u_k:=a_k\circ h.
\]
The convolution is defined even though $H_m$ is not integrable
on the whole line: at each $t$ it integrates only over a compact
interval, where $H_m$ is integrable. Differentiation under the
integral is valid on a common compact set as $t$ varies in a
bounded neighborhood, because each derivative of the kernel
is bounded and compactly supported. Thus $a_k$ is smooth.
Applying the distributional identity to the test function
$s\mapsto\rho_{\varepsilon_k}(t-s)$, the signs from differentiating
the kernel cancel, giving
\[
 a_k^{(j)}(t)
 =\int_{\mathbb R}\rho_{\varepsilon_k}(t-s)H_{m-j}(s)\,ds,
 \qquad 0\leq j\leq m.
\]
Consequently, $u_k$ is a global smooth function on $S^2$.
Its derivatives are continuous and the sphere is compact with
finite volume, so $u_k\in C^{m,p}(S^2)$.

For $r\geq1$, $H_r$ is Lipschitz on $[-2,2]$ with constant
\[
 L_r:=\frac{2^{r-1}}{(r-1)!}.
\]
On each half-line this follows from the fundamental theorem
of calculus and $|H_r'|\leq L_r$; if a segment crosses zero,
split it in two and add the estimates. For $|t|\leq1$ and
$0<\varepsilon<1$, every argument $t-s$ of the kernel lies
in $[-2,2]$, and hence
\[
 |(\rho_\varepsilon*H_r)(t)-H_r(t)|
 \leq\int_{\mathbb R}\rho_\varepsilon(s)L_r|s|\,ds
 \leq L_r\varepsilon.
\]
For $r=0$, the convolution takes values between zero and one,
and equals $H_0(t)$ when $|t|\geq\varepsilon$.
The difference has absolute value at most one and is supported
in $(-\varepsilon,\varepsilon)$.

The coordinates $(\theta,z)$ and the volume element calculated
in Example~\ref{ej:sobolev-altura-y-esquina-esfera} give, for
any nonnegative measurable function $b$ on $[-1,1]$,
\[
 \int_{S^2}b(h)\,d\lambda_{\mathbf g}
 =2\pi\int_{-1}^1b(t)\,dt.
\]
In particular, $\lambda_{\mathbf g}(S^2)=4\pi$, and the band
$|h|<\varepsilon$ has area $4\pi\varepsilon$. The preceding
estimates imply, for $0\leq j\leq m$,
\[
 \|a_k^{(j)}(h)-H_{m-j}(h)\|_{L^p(S^2)}
 \leq
 \begin{cases}
 (4\pi)^{\frac{1}{p}}L_{m-j}\varepsilon_k,&j<m,\\
 (4\pi\varepsilon_k)^{\frac{1}{p}},&j=m.
 \end{cases}
\]
The second line includes $m=0$.

There are smooth tensors $B_{s,j}$ of type $(0,s)$, defined
for $0\leq j\leq s$, such that for every smooth function $a$,
\[
 \nabla^s(a\circ h)=\displaystyle\sum_{j=0}^s a^{(j)}(h)B_{s,j}.
\]
At order zero take $B_{0,0}=1$. If the formula holds at order
$s$, the product rule and the new slot in first position give
\[
 B_{s+1,j}=\nabla B_{s,j}+dh\otimes B_{s,j-1},
\]
where terms with indices out of range mean zero. This recurrence
proves the formula by induction and the smoothness of every tensor.
Set
\[
 T_s:=\displaystyle\sum_{j=0}^s H_{m-j}(h)B_{s,j},\qquad
 A_{s,j}:=\displaystyle\sup_{S^2}|B_{s,j}|_{\mathbf g},\qquad 0\leq s\leq m.
\]
The constants are finite by compactness. Each $T_s$ is measurable
and bounded, hence belongs to $L^p(S^2,T^{(0,s)}(TS^2))$.
The triangle inequality gives
\[
 \|\nabla^su_k-T_s\|_{L^p(S^2)}
 \leq\displaystyle\sum_{j=0}^s A_{s,j}
       \|a_k^{(j)}(h)-H_{m-j}(h)\|_{L^p(S^2)}\longrightarrow0.
\]
There are finitely many orders. If
$c_k:=\displaystyle\sum_{s=0}^m\|\nabla^su_k-T_s\|_{L^p(S^2)}$, then
$c_k\to0$, and
\[
 \|u_k-u_\ell\|_{W^{m,p}(S^2)}
 \leq\displaystyle\sum_{s=0}^m\|\nabla^su_k-\nabla^su_\ell\|_{L^p(S^2)}
 \leq c_k+c_\ell.
\]
This proves the Cauchy property. In the completion its limit
is the class of this sequence; the continuous order-zero map into
$L^p$ sends that class to $u:=H_m\circ h$, by the order-zero
convergence just proved. Continuity of this map follows
from $\|f\|_{L^p}\leq\|f\|_{W^{m,p}}$ and the universal
property of completion.

If $m\geq1$, $u$ is continuous. A smooth function $w$ that were
a limit in the Sobolev norm would equal $u$ almost everywhere,
hence everywhere by continuity: every nonempty open subset
of $S^2$ has positive measure, since its volume density in a
chart is positive. Consider the smooth curve
\[
 \gamma(t)=(\sqrt{1-t^2},0,t),\qquad |t|<1/2.
\]
We have $u(\gamma(t))=H_m(t)$. For $t<0$ its derivative of
order $m$ is zero, and for $t>0$ it is one. Thus $H_m$ is not
of class $C^m$ at zero and cannot be $w\circ\gamma$.
If $m=0$, equality almost everywhere would force $w=0$ on
the hemisphere $h<0$ and $w=1$ on the hemisphere $h>0$, by
the same continuity and positive-measure argument on each
hemisphere. The limits of $w(\gamma(t))$ as $t$ tends to zero
from both sides would differ, contradicting continuity of $w$.
In every case $C^{m,p}(S^2)$ is incomplete in the Sobolev norm,
and the limit added by completion has no smooth representative.
\end{example}

\begin{lemma}[Localizing Euclidean examples]
\label{lem:localizar-ejemplos-sobolev-escalar}
Let $(M,\mathbf g)$ be a Riemannian manifold without boundary of dimension
$n\geq1$,
$\phi\colon U\to V\subseteq\mathbb R^n$ a chart, and
$1\leq p<\infty$, $f\in W^{1,p}(V)$, $\chi\in C_c^\infty(V)$.
Then the function
\[
 u(p_0):=
 \begin{cases}
  (\chi f)(\phi(p_0)),&p_0\in U,\\
  0,&p_0\notin U
 \end{cases}
\]
represents a compactly supported element of $W^{1,p}(M)$.
Conversely, if $u\in W^{1,p}(M)$ and $V'\Subset V$ is open, its expression
$u\circ\phi^{-1}$ belongs to $W^{1,p}(V')$.
The local condition holds in every chart, but alone does not replace
the global integrability required by the norm.
\end{lemma}
\begin{proof}
By the weak product rule and zero extension in Chapter~3,
$w:=\chi f$, extended by zero, belongs to $W^{1,p}(\mathbb R^n)$ and
\[
 D_iw=(D_i\chi)f+\chi D_if.
\]
Its support $K$ is compact and contained in $V$.
Choose $\delta_0>0$ such that the compact set
$K':=\{x\mid\operatorname{dist}_{\mathrm{euc}}(x,K)\leq\delta_0\}$
is contained in $V$; if $K$ is empty, the conclusion is immediate.
The mollifiers of Theorem~\ref{teo:molificadores-euclidianos} give smooth
$w_\delta:=w*\rho_\delta$ supported in $K'$ for $0<\delta<\delta_0$.
The identity
\[
 D_iw_\delta=(D_iw)*\rho_\delta
\]
follows by applying the weak-derivative definition to
$y\mapsto\rho_\delta(x-y)$, a test function.
Approximation of the identity applied to $w$ and to each $D_iw$ gives
$w_\delta\to w$ in $W^{1,p}(\mathbb R^n)$.
Local convergence can be applied on a compact set containing all
supports, since the functions vanish outside it.

Transport $w_\delta$ through $\phi$ and extend by zero to obtain
$u_\delta\in C_c^\infty(M)$.
On $K'$, let $\rho(x):=\sqrt{\det(g_{ij}(x))}\leq A$ and
$g^{ij}(x)\xi_i\xi_j\leq B|\xi|_2^2$, with $A,B>0$.
These constants exist by compactness and positivity.
For any smooth function $z$ supported in $K'$, its transport
$\widetilde z$ satisfies
\[
 \begin{aligned}
 \|\widetilde z\|_{W^{1,p}(M)}^p
 &\leq A\left(\|z\|_{L^p(\mathbb R^n)}^p+
 B^{p/2}n^{\max\{0,p/2-1\}}
 \sum_{i=1}^n\|D_i z\|_{L^p(\mathbb R^n)}^p\right)\\
 &\leq C^p\|z\|_{W^{1,p}(\mathbb R^n)}^p,
 \end{aligned}
\]
where $C^p:=A\max\{1,B^{p/2}n^{\max\{0,p/2-1\}}\}$.
Here we used the finite-dimensional inequality
$|\xi|_2\leq n^{\max\{0,1/2-1/p\}}
(\displaystyle\sum_{i=1}^n|\xi_i|^p)^{1/p}$.
Applying the estimate to $w_\delta-w_{\delta'}$ shows that $(u_\delta)$
is Cauchy in the completion and that its $L^p$ limit is the stated $u$.
The differentials converge to the transport of the weak derivatives
of $w$, so the class is well identified.

For the converse, approximate $u$ in the completion norm by smooth
functions. On $\overline{V'}$, the density $\rho$ has a positive lower
bound, and the matrix $(g^{ij})$ has eigenvalues uniformly bounded
above and away from zero. If $\rho\geq a>0$ and
$g^{ij}\xi_i\xi_j\geq b|\xi|_2^2$ there, a smooth function $z$ satisfies
the explicit estimate
\[
 \|z\circ\phi^{-1}\|_{W^{1,p}(V')}^p
 \leq a^{-1}\max\{1,b^{-p/2}n^{\max\{0,1-p/2\}}\}
     \|z\|_{W^{1,p}(M)}^p.
\]
Here we use
$\displaystyle\sum_{i=1}^n|\xi_i|^p
\leq n^{\max\{0,1-p/2\}}|\xi|_2^p$.
The representations of the approximating sequence are therefore Cauchy
in $W^{1,p}(V')$. Completeness proved in Chapter~3 and convergence in
$L^p(V')$ identify their limit with $u\circ\phi^{-1}$.

Finally, on a noncompact manifold a function may be Sobolev on every
relatively compact region while having infinite global norm.
For example, $1$ is smooth on $\mathbb R^n$ and belongs to
$W^{1,p}(V')$ for every relatively compact open set $V'$, but does
not belong to $L^p(\mathbb R^n)$ for any $p<\infty$.
\end{proof}

\begin{example}[A corner, a singularity, and a jump on the manifold]
\label{ej:sobolev-localizados-variedad}
In a chart centered at $p_0$, choose a coordinate ball
$B_{\mathrm{euc}}(0,R)\subseteq V$ and a function
$\chi\in C_c^\infty(B_{\mathrm{euc}}(0,R))$ equal to $1$ near $0$.
If $n\geq1$, the local formula $\chi(x)|x^1|$, extended by zero on $M$,
is continuous and belongs to $W^{1,p}(M)$ for every $1\leq p<\infty$.
The weak derivative of $|x^1|$ is $\operatorname{sgn}(x^1)$ in
direction $1$ and zero in all other directions, as follows from
Example~\ref{ej:sobolev-potencias-intervalo} using Fubini.
The product rule and the preceding lemma give global membership.
It is not $C^1$ where $\chi=1$ and $x^1=0$.

If $n=3$, the local formula $\chi(x)|x|^{-1/4}$ produces an element
of $W^{1,2}(M)$ with no continuous representative near $p_0$.
Example~\ref{ej:sobolev-singularidad-radial} applies after a dilation
to the ball of radius $R$, and the lemma transports it.
The Riemannian density is positive and bounded above and away from zero
on a smaller ball; hence sets of positive Euclidean measure there have
positive Riemannian measure. The singularity cannot be eliminated by
changing values on a null set.

Finally, on $S^2$, the function $b=\mathbf 1_{\{h>0\}}$ of
Example~\ref{ej:regularidad-escalar-esfera} does not belong to
$W^{1,p}(S^2)$ for any $1\leq p<\infty$.
Near the equator, use the coordinates $(\theta,z)$ from the previous
example, restricting $\theta$ to an interval to obtain a chart.
If $b$ belonged to the space, the converse of the lemma would give a
locally integrable weak derivative $D_zb$.
Using a test function $\eta(\theta)\varphi(z)$, with
$\eta\in C_c^\infty$ and $\displaystyle\int\eta=1$, Fubini would give
\[
 \int H(z)\varphi'(z)\,dz
 =-\int\left(\int D_zb(\theta,z)\eta(\theta)\,d\theta\right)\varphi(z)\,dz.
\]
The function in parentheses belongs to $L^1_{\mathrm{loc}}$ by Fubini,
contradicting Example~\ref{ej:sobolev-salto-y-representantes}.
Charts recognize the corner and jump; the metric measures integrability,
and compatibility ensures that the objects are global.
\end{example}

\section{Kato's inequality}
We shall prove Rademacher's theorem on Riemannian manifolds in order to establish an important inequality known as \textit{Kato's inequality} for sections of vector bundles. It will be useful in proving the Sobolev embeddings in the next section.

Recall Rademacher's theorem in Euclidean spaces.
\begin{theorem}[Rademacher in Euclidean spaces]\label{rademacher euclidiano}
Let $U\subseteq\mathbb{R}^{n}$ be open and let $f\colon U\longrightarrow\mathbb{R}^{m}$ be a Lipschitz function. Then $f$ is differentiable almost everywhere on $U$.
\end{theorem}
This theorem is proved in \cite{Rademacher} and allows us to prove a version of Kato's inequality. We shall establish its counterpart on Riemannian manifolds to develop nonlinear analysis on Riemannian manifolds further.
\begin{theorem}[Rademacher on Riemannian manifolds]\label{rademacher en variedades}
Let $(M,\mathbf{g})$ be a Riemannian manifold with or without boundary and let $u\colon M\longrightarrow\mathbb{R}$ be a Lipschitz function. Then the first covariant derivative $\nabla u$ exists almost everywhere on $M$.
\end{theorem}
\begin{proof}
By Proposition \ref{bolacoordenadaregular}, $M$ has a countable basis of regular coordinate balls (and if $M$ has boundary, some are regular coordinate half-balls). By the definition of regular coordinate balls (or half-balls), we may assume that there is a countable atlas $\{(B_{j},\phi_{j})\mid j\in\mathbb{N}\}$ of regular coordinate balls (or half-balls) such that $\{(\overline{B}_{j},\phi_{j})\mid j\in\mathbb{N}\}$ is a family of compact regular coordinate balls (or half-balls), and $M=\displaystyle\bigcup_{j\in\mathbb{N}}B_{j}=\displaystyle\bigcup_{j\in\mathbb{N}}\overline{B}_{j}$. Each $\overline{B}_{j}$ is a compact subset contained in an ambient chart; it need not be viewed as a smooth submanifold with boundary. For each $j\in\mathbb{N}$, let $C_{j}=\phi_{j}(\overline{B}_{j})$.\\

The set $B_{j}$ admits two Riemannian metrics, $\mathbf{g}$ and $\boldsymbol{\delta}:=\phi_{j}^{*}\overline{\mathbf{g}}$, where $\overline{\mathbf{g}}$ denotes the usual metric on $\mathbb{R}^{n}$. By the pointwise comparison of Lemma~\ref{lema cg} on the compact set $\overline{B_j}$, in the sense of Remark~\ref{obs:comparacion-metricas-compactos-coordenados}, there is $K_j>0$ such that for any $p\in B_{j}$ and $v\in T_{p}M$, $\mathbf{g}_{p}(v,v)\leq K_{j}\boldsymbol{\delta}_{p}(v,v)$.\\

We shall prove that $h_{j}:=u\circ\phi_{j}^{-1}\colon C_{j}\longrightarrow\mathbb{R}$ is Lipschitz. Since $u$ is Lipschitz, there is $K>0$ such that $|u(p)-u(q)|\leq Kd_{\mathbf{g}}(p,q)$ for any $p,q\in M$. Let $x=\phi_{j}(p),y=\phi_{j}(q)\in C_{j}$. Since $C_{j}$ is a closed ball (or closed half-ball) in $\mathbb{R}^{n}$, it is convex, and hence the segment joining $x$ to $y$ remains in $C_{j}$. Let $\gamma\colon [0,1]\longrightarrow C_{j}$ be given by $\gamma(t)=(1-t)x+ty$. Then $\gamma$ is a smooth curve with $\gamma(0)=x$ and $\gamma(1)=y$. Thus $\widetilde{\gamma}_{j}:=\phi_{j}^{-1}\circ\gamma\colon [0,1]\longrightarrow \overline{B}_{j}$ is a smooth curve with $\widetilde{\gamma}_{j}(0)=p$ and $\widetilde{\gamma}_{j}(1)=q$.\\

We therefore have
$$|h_{j}(x)-h_{j}(y)|=|u(\phi_{j}^{-1}(x))-u(\phi_{j}^{-1}(y))|\leq Kd_{\mathbf{g}}(p,q)\leq KL_{\mathbf{g}}(\widetilde{\gamma}_{j})$$
$$=K\int_{0}^{1}|\widetilde{\gamma}_{j}'(t)|_{\mathbf{g}}dt=K\int_{0}^{1}\sqrt{\mathbf{g}_{\widetilde{\gamma}_{j}(t)}(\widetilde{\gamma}_{j}'(t),\widetilde{\gamma}_{j}'(t))}dt$$
$$\leq K\sqrt{K_{j}}\int_{0}^{1}\sqrt{\boldsymbol{\delta}_{(\phi_{j}^{-1}\circ\gamma)(t)}((\phi_{j}^{-1}\circ\gamma)'(t),(\phi_{j}^{-1}\circ\gamma)'(t))}dt$$
$$=K\sqrt{K_{j}}\int_{0}^{1}\sqrt{\overline{\mathbf{g}}_{\gamma(t)}(\gamma'(t),\gamma'(t))}dt=K\sqrt{K_{j}}\int_{0}^{1}\|\gamma'(t)\|dt$$
$$=K\sqrt{K_{j}}\int_{0}^{1}\|y-x\|dt=K\sqrt{K_{j}}\|x-y\|,$$
which proves that $h_{j}$ is Lipschitz. By Rademacher's theorem (Theorem \ref{rademacher euclidiano}), the derivative of $h_{j}$ exists almost everywhere on $\operatorname{int}(C_{j})$. But $C_j$ is a closed ball or a closed half-ball, so its Euclidean boundary has zero $n$-dimensional measure. Thus the derivative of $h_{j}$ exists almost everywhere on $C_{j}$; that is, there exists $N_{j}\subseteq C_{j}$ such that $\lambda_{n}(N_{j})=0$ and $dh_{j}$ exists on $C_{j}\setminus N_{j}$. Let $N_{j}'=\phi_{j}^{-1}(N_{j})$. Since
$$\lambda_{\mathbf{g}}(N_{j}')=\int_{\phi_{j}(N_{j}')}(\phi_{j})_{*}\sqrt{\det(\mathbf{g})}d\lambda_{n}=\int_{N_{j}}\sqrt{\det(\mathbf{g})}\circ\phi_{j}^{-1}d\lambda_{n}=0,$$
we conclude that $N_{j}'$ has zero Riemann--Lebesgue measure in $M$.\\

Let $p\in B_{j}\setminus N_{j}'$. Then $\phi_{j}(p)\in C_{j}\setminus N_{j}$, so $dh_{j}|_{\phi_{j}(p)}$ exists. Since $\phi_{j}$ is smooth, $d\phi_{j}|_{p}$ exists. Because $u|_{B_{j}}=h_{j}\circ\phi_{j}$, the chain rule shows that $du|_{p}$ exists. Since the first covariant derivative of a function agrees with its differential, $\nabla u(p)$ exists. Thus $\nabla u$ exists on $B_{j}\setminus N_{j}'$. But $N=\displaystyle\bigcup_{j\in\mathbb{N}}N_{j}'$ satisfies $\lambda_{\mathbf{g}}(N)=0$, and moreover $\nabla u(p)$ exists for every $p\in\displaystyle\bigcup_{j\in\mathbb{N}}B_{j}\setminus\displaystyle\bigcup_{j\in\mathbb{N}}N_{j}'=M\setminus N$. Therefore, $\nabla u$ exists almost everywhere on $M$.
\end{proof}

We now show that if $\mathbf{u}\in\Gamma(\mathbf{E})$, then $\nabla|\nabla^{s}\mathbf{u}|_{\mathbf{g},\mathbf{h}_{\mathbf{E}}}$ exists almost everywhere on $M$. First, we prove two results.

\begin{proposition}\label{D^s es lipschitz}
Let $U\subseteq\mathbb{R}^{n}$ be a convex open subset, let $s\in\mathbb{N}\cup\{0\}$, and let $f\colon U\longrightarrow\mathbb{R}$ be a function of class $C^{s+1}$. If there is $M>0$ such that $\|\nabla^{s+1}f(x)\|_{F}\leq M$ for every $x\in U$, then
\[
    \|\nabla^{s}f(x)-\nabla^{s}f(y)\|_{F}\leq M\|x-y\|
    \]
for any $x,y\in U$, where $\nabla^{0}f=f$ and, for $s\geq1$,
\[
    \nabla^{s}f(x)=\displaystyle\frac{\partial^{s}f}{\partial x_{i_{1}}\cdots\partial x_{i_{s}}}\Biggr|_{x}\mathbf{d}x^{i_{1}}\otimes\cdots\otimes \mathbf{d}x^{i_{s}}.
    \]
\end{proposition}
\begin{proof}
If $s=0$, the following argument is interpreted with $\nabla^{0}f=f$. Let $x,y\in U$. Define $\gamma\colon [0,1]\longrightarrow U$ by $\gamma(t)=x+t(y-x)$. Then $\gamma'(t)=y-x$. We have
$$\nabla^{s}f(y)-\nabla^{s}f(x)=\displaystyle\frac{\partial^{s}f}{\partial x_{i_{1}}\cdots\partial x_{i_{s}}}\Biggr|_{\gamma(1)}\mathbf{d}x^{i_{1}}\otimes\cdots\otimes \mathbf{d}x^{i_{s}}-\displaystyle\frac{\partial^{s}f}{\partial x_{i_{1}}\cdots\partial x_{i_{s}}}\Biggr|_{\gamma(0)}\mathbf{d}x^{i_{1}}\otimes\cdots\otimes \mathbf{d}x^{i_{s}}$$
$$=\left(\displaystyle\frac{\partial^{s}f}{\partial x_{i_{1}}\cdots\partial x_{i_{s}}}(\gamma(1))-\displaystyle\frac{\partial^{s}f}{\partial x_{i_{1}}\cdots\partial x_{i_{s}}}(\gamma(0))\right)\mathbf{d}x^{i_{1}}\otimes\cdots\otimes \mathbf{d}x^{i_{s}}$$
$$=\left(\int_{0}^{1}\displaystyle\frac{d}{dt}\left(\displaystyle\frac{\partial^{s}f}{\partial x_{i_{1}}\cdots\partial x_{i_{s}}}(\gamma(t))\right)dt\right)\mathbf{d}x^{i_{1}}\otimes\cdots\otimes \mathbf{d}x^{i_{s}}$$
$$=\left(\int_{0}^{1}\left\langle \nabla\left(\displaystyle\frac{\partial^{s}f}{\partial x_{i_{1}}\cdots\partial x_{i_{s}}}\right)_{\gamma(t)},y-x\right\rangle_{\mathbb{R}^{n}}dt\right)\mathbf{d}x^{i_{1}}\otimes\cdots\otimes \mathbf{d}x^{i_{s}}$$
$$=\left(\int_{0}^{1}\displaystyle\sum_{j=1}^{n}\displaystyle\frac{\partial^{s+1}f}{\partial x_{j}\partial x_{i_{1}}\cdots\partial x_{i_{s}}}(\gamma(t))(y^{j}-x^{j})dt\right)\mathbf{d}x^{i_{1}}\otimes\cdots\otimes \mathbf{d}x^{i_{s}}.$$
Thus,
$$\|\nabla^{s}f(y)-\nabla^{s}f(x)\|_{F}=\sqrt{\displaystyle\sum_{i_{1},\dots,i_{s}=1}^{n}\left(\int_{0}^{1}\displaystyle\sum_{j=1}^{n}\displaystyle\frac{\partial^{s+1}f}{\partial x_{j}\partial x_{i_{1}}\cdots\partial x_{i_{s}}}(\gamma(t))(y^{j}-x^{j})dt\right)^{2}}$$
$$\leq\sqrt{\displaystyle\sum_{i_{1},\dots,i_{s}=1}^{n}\int_{0}^{1}\left(\displaystyle\sum_{j=1}^{n}\displaystyle\frac{\partial^{s+1}f}{\partial x_{j}\partial x_{i_{1}}\cdots\partial x_{i_{s}}}(\gamma(t))(y^{j}-x^{j})\right)^{2}dt},$$
where the last inequality follows from Jensen's inequality (Proposition \ref{desigualdad de Jensen}) and convexity of $\varphi(x)=x^{2}$. But
$$\sqrt{\displaystyle\sum_{i_{1},\dots,i_{s}=1}^{n}\int_{0}^{1}\left(\displaystyle\sum_{j=1}^{n}\displaystyle\frac{\partial^{s+1}f}{\partial x_{j}\partial x_{i_{1}}\cdots\partial x_{i_{s}}}(\gamma(t))(y^{j}-x^{j})\right)^{2}dt}$$
$$=\sqrt{\displaystyle\sum_{i_{1},\dots,i_{s}=1}^{n}\int_{0}^{1}\left\langle\left(\displaystyle\frac{\partial^{s+1}f}{\partial x_{1}\partial x_{i_{1}}\cdots\partial x_{i_{s}}}(\gamma(t)),\dots,\displaystyle\frac{\partial^{s+1}f}{\partial x_{n}\partial x_{i_{1}}\cdots\partial x_{i_{s}}}(\gamma(t))\right),y-x\right\rangle_{\mathbb{R}^{n}}^{2}dt}$$
$$\leq\sqrt{\displaystyle\sum_{i_{1},\dots,i_{s}=1}^{n}\int_{0}^{1}\left\|\left(\displaystyle\frac{\partial^{s+1}f}{\partial x_{1}\partial x_{i_{1}}\cdots\partial x_{i_{s}}}(\gamma(t)),\dots,\displaystyle\frac{\partial^{s+1}f}{\partial x_{n}\partial x_{i_{1}}\cdots\partial x_{i_{s}}}(\gamma(t))\right)\right\|^{2}\|y-x\|^{2}dt}$$
$$=\|x-y\|\sqrt{\displaystyle\sum_{i_{1},\dots,i_{s}=1}^{n}\int_{0}^{1}\left\|\left(\displaystyle\frac{\partial^{s+1}f}{\partial x_{1}\partial x_{i_{1}}\cdots\partial x_{i_{s}}}(\gamma(t)),\dots,\displaystyle\frac{\partial^{s+1}f}{\partial x_{n}\partial x_{i_{1}}\cdots\partial x_{i_{s}}}(\gamma(t))\right)\right\|^{2}dt}$$
$$=\|x-y\|\sqrt{\int_{0}^{1}\displaystyle\sum_{i_{1},\dots,i_{s}=1}^{n}\left\|\left(\displaystyle\frac{\partial^{s+1}f}{\partial x_{1}\partial x_{i_{1}}\cdots\partial x_{i_{s}}}(\gamma(t)),\dots,\displaystyle\frac{\partial^{s+1}f}{\partial x_{n}\partial x_{i_{1}}\cdots\partial x_{i_{s}}}(\gamma(t))\right)\right\|^{2}dt}$$
$$=\|x-y\|\sqrt{\int_{0}^{1}\displaystyle\sum_{i_{1},\dots,i_{s}=1}^{n}\displaystyle\sum_{j=1}^{n}\left(\displaystyle\frac{\partial^{s+1}f}{\partial x_{j}\partial x_{i_{1}}\cdots\partial x_{i_{s}}}(\gamma(t))\right)^{2}dt}$$
$$=\|x-y\|\sqrt{\int_{0}^{1}\|\nabla^{s+1}f(\gamma(t))\|_{F}^{2}dt}\leq\|x-y\|\sqrt{\int_{0}^{1}M^{2}dt}=M\|x-y\|.$$
\end{proof}

\begin{lemma}\label{lema para ver que existe casi en todas partes}
Let $U\subseteq\mathbb{R}^{n}$ be open with $\overline{B}_{\mathrm{euc}}(0,2)\subseteq U$, and let $f\colon U\longrightarrow[0,\infty)$ be a function of class $C^{2}$. Then the function $\sqrt{f}$ is Lipschitz on $\overline{B}_{\mathrm{euc}}(0,1)$.
\end{lemma}
\begin{proof}
Since $f$ is of class $C^{2}$ on the compact set $\overline{B}_{\mathrm{euc}}(0,2)$, its gradient and Hessian matrix are bounded on that set. Let
$$K:=1+\sup_{z\in\overline{B}_{\mathrm{euc}}(0,2)}\|\nabla f(z)\|+\sup_{z\in\overline{B}_{\mathrm{euc}}(0,2)}\|\nabla^{2}f(z)\|_{F}.$$
Then $K>0$ and $\|\nabla f(z)\|\leq K$, $\|\nabla^{2}f(z)\|_{F}\leq K$ for every $z\in\overline{B}_{\mathrm{euc}}(0,2)$.

For $\varepsilon>0$, consider $f_{\varepsilon}=f+\varepsilon$. Fix $x\in\overline{B}_{\mathrm{euc}}(0,1)$. If $\nabla f(x)=0$, then $\nabla\sqrt{f_{\varepsilon}}(x)=0$. Thus assume $\nabla f(x)\neq0$ and define
$$v:=\displaystyle\frac{\nabla f(x)}{\|\nabla f(x)\|}.$$
Then $\|v\|=1$ and, since $x\in\overline{B}_{\mathrm{euc}}(0,1)$, we have $x-tv\in\overline{B}_{\mathrm{euc}}(0,2)$ for every $t\in[0,1]$.

Consider the function $g\colon [0,1]\longrightarrow(0,\infty)$ given by $g(t):=f_{\varepsilon}(x-tv)$. Then
$$g'(0)=-\langle\nabla f(x),v\rangle_{\mathbb{R}^{n}}=-\|\nabla f(x)\|.$$
Moreover, for every $t\in[0,1]$,
$$g''(t)=\displaystyle\sum_{i,j=1}^{n}\displaystyle\frac{\partial^{2}f}{\partial x_{i}\partial x_{j}}(x-tv)v_{i}v_{j}.$$
The Cauchy--Schwarz inequality from Theorem~\ref{teo:cauchy-schwarz-hilbert}, applied in $\mathbb{R}^{n^{2}}$, gives
$$|g''(t)|\leq\left(\displaystyle\sum_{i,j=1}^{n}\left|\displaystyle\frac{\partial^{2}f}{\partial x_{i}\partial x_{j}}(x-tv)\right|^{2}\right)^{\frac{1}{2}}\left(\displaystyle\sum_{i,j=1}^{n}|v_{i}v_{j}|^{2}\right)^{\frac{1}{2}}$$
$$=\|\nabla^{2}f(x-tv)\|_{F}\|v\|^{2}\leq K.$$
By Taylor's formula with integral remainder from Theorem~\ref{taylor multivariable multiindices}, applied in one dimension, for every $t\in[0,1]$,
$$g(t)=g(0)+g'(0)t+\int_{0}^{t}(t-r)g''(r)dr\leq g(0)-\|\nabla f(x)\|t+\displaystyle\frac{K}{2}t^{2}.$$

Let $t_{0}:=\frac{\|\nabla f(x)\|}{K}$. Since $\|\nabla f(x)\|\leq K$, we have $t_{0}\in[0,1]$. Moreover, since $f\geq0$ and $\varepsilon>0$, we have $g(t_{0})=f(x-t_{0}v)+\varepsilon\geq\varepsilon>0$. Therefore,
$$0<g(t_{0})\leq g(0)-\displaystyle\frac{\|\nabla f(x)\|^{2}}{K}+\displaystyle\frac{K}{2}\displaystyle\frac{\|\nabla f(x)\|^{2}}{K^{2}}=g(0)-\displaystyle\frac{\|\nabla f(x)\|^{2}}{2K}.$$
It follows that $\|\nabla f(x)\|^{2}<2Kg(0)=2K(f(x)+\varepsilon)$ and consequently
$$\|\nabla\sqrt{f_{\varepsilon}}(x)\|=\displaystyle\frac{\|\nabla f(x)\|}{2\sqrt{f(x)+\varepsilon}}\leq\sqrt{\displaystyle\frac{K}{2}}.$$

This bound holds for every $x\in\overline{B}_{\mathrm{euc}}(0,1)$, and the constant on the right-hand side is independent of $x$ and $\varepsilon$. Since $\sqrt{f_{\varepsilon}}$ is of class $C^{1}$ on the convex open subset $B_{\mathrm{euc}}(0,1)$, Proposition \ref{D^s es lipschitz}, applied with $s=0$, implies
$$\left|\sqrt{f(x)+\varepsilon}-\sqrt{f(y)+\varepsilon}\right|\leq\sqrt{\displaystyle\frac{K}{2}}\|x-y\|$$
for any $x,y\in B_{\mathrm{euc}}(0,1)$. By continuity, the inequality also holds for any $x,y\in\overline{B}_{\mathrm{euc}}(0,1)$. Finally, letting $\varepsilon\to0$, we obtain
$$|\sqrt{f(x)}-\sqrt{f(y)}|\leq\sqrt{\displaystyle\frac{K}{2}}\|x-y\|$$
for any $x,y\in\overline{B}_{\mathrm{euc}}(0,1)$, so $\sqrt{f}$ is Lipschitz on $\overline{B}_{\mathrm{euc}}(0,1)$.
\end{proof}

\begin{lemma}\label{lema para kato}
Let $(M,\mathbf{g})$ be a Riemannian manifold with or without boundary, and let $\mathbf{E}\to M$ be a smooth real vector bundle equipped with a bundle metric $\mathbf{h}_{\mathbf{E}}$ and a compatible connection $\nabla^{\mathbf{E}}$. If $\mathbf{u}\in\Gamma(\mathbf{E})$, then for every integer $s\geq0$, the function $|\nabla^s\mathbf{u}|_{\mathbf{g},\mathbf{h}_{\mathbf{E}}}$ is locally Lipschitz on $\operatorname{Int}(M)$. In particular, $\nabla|\nabla^s\mathbf{u}|_{\mathbf{g},\mathbf{h}_{\mathbf{E}}}$ exists almost everywhere on $M$.
\end{lemma}

\begin{proof}
Since $\partial M$ has measure zero, it suffices to prove that $\nabla|\nabla^s\mathbf{u}|_{\mathbf{g},\mathbf{h}_{\mathbf{E}}}$ exists almost everywhere on $\operatorname{Int}(M)$.

By Lemma~\ref{lem:vecindades-normales-convexas-locales}, choose
a strongly convex regular normal geodesic coordinate ball $(V,\phi)$ whose closure is contained in $\operatorname{Int}(M)$. Composing $\phi$ with an affine transformation of $\mathbb R^n$, we may assume $\overline{B}_{\mathrm{euc}}(0,3)\subseteq\phi(V)$. Define $B:=\phi^{-1}(\overline{B}_{\mathrm{euc}}(0,1))$.

On $V$, we have
\[
|\nabla^s\mathbf{u}|_{\mathbf{g},\mathbf{h}_{\mathbf{E}}}
=
\sqrt{g^{i_1j_1}\cdots g^{i_sj_s}h_{ab}
(\nabla^s\mathbf{u})^a_{i_1\cdots i_s}
(\nabla^s\mathbf{u})^b_{j_1\cdots j_s}}
=
\sqrt{\langle\nabla^s\mathbf{u},\nabla^s\mathbf{u}\rangle_{\mathbf{g},\mathbf{h}_{\mathbf{E}}}},
\]
where, if $s=0$, the factors corresponding to $\mathbf{g}^{-1}$ are absent. Since $\langle\nabla^s\mathbf{u},\nabla^s\mathbf{u}\rangle_{\mathbf{g},\mathbf{h}_{\mathbf{E}}}\circ\phi^{-1}$ is smooth and $\phi(V)$ contains $\overline{B}_{\mathrm{euc}}(0,2)$, Lemma \ref{lema para ver que existe casi en todas partes} implies that
\[
|\nabla^s\mathbf{u}|_{\mathbf{g},\mathbf{h}_{\mathbf{E}}}\circ\phi^{-1}
=
\sqrt{\langle\nabla^s\mathbf{u},\nabla^s\mathbf{u}\rangle_{\mathbf{g},\mathbf{h}_{\mathbf{E}}}\circ\phi^{-1}}
\]
is Lipschitz on $\overline{B}_{\mathrm{euc}}(0,1)$. Thus there is $C>0$ such that
\[
\left|
\bigl(|\nabla^s\mathbf{u}|_{\mathbf{g},\mathbf{h}_{\mathbf{E}}}\circ\phi^{-1}\bigr)(x)
-
\bigl(|\nabla^s\mathbf{u}|_{\mathbf{g},\mathbf{h}_{\mathbf{E}}}\circ\phi^{-1}\bigr)(y)
\right|
\leq
C\|x-y\|
\]
for any $x,y\in\overline{B}_{\mathrm{euc}}(0,1)$.

Let $p,q\in B$. Let $\boldsymbol{\delta}:=\phi^{*}\overline{\mathbf{g}}$ be the Euclidean metric transported to $V$. Since $\overline V$ is compact, Lemma~\ref{lema cg}, applied to the metrics $\mathbf{g}$ and $\boldsymbol{\delta}$, provides a constant $K>0$ such that $|X|_{\boldsymbol{\delta}}\leq K|X|_{\mathbf{g}}$ for every $X\in TV$. Since $V$ is strongly convex, there is a minimizing geodesic $\gamma\colon [0,1]\longrightarrow V$ joining $p$ to $q$ and satisfying $L_{\mathbf{g}}(\gamma)=d_{\mathbf{g}}(p,q)$. The curve $\phi\circ\gamma$ joins $\phi(p)$ to $\phi(q)$, and therefore
\[
\|\phi(p)-\phi(q)\|
\leq
L_{\overline{\mathbf{g}}}(\phi\circ\gamma)
=
L_{\boldsymbol{\delta}}(\gamma)
\leq
K L_{\mathbf{g}}(\gamma)
=
K d_{\mathbf{g}}(p,q).
\] Consequently,
\[
\begin{aligned}
\left|
|\nabla^s\mathbf{u}|_{\mathbf{g},\mathbf{h}_{\mathbf{E}}}(p)-|\nabla^s\mathbf{u}|_{\mathbf{g},\mathbf{h}_{\mathbf{E}}}(q)
\right|
&\leq
C\|\phi(p)-\phi(q)\|\\
&\leq
CK\,d_{\mathbf{g}}(p,q),
\end{aligned}
\]
so $|\nabla^s\mathbf{u}|_{\mathbf{g},\mathbf{h}_{\mathbf{E}}}$ is Lipschitz on $B$.

By Lemma~\ref{lem:vecindades-normales-convexas-locales}, the interior
admits a cover by such strongly convex normal coordinate balls, we conclude that $|\nabla^s\mathbf{u}|_{\mathbf{g},\mathbf{h}_{\mathbf{E}}}$ is locally Lipschitz on $\operatorname{Int}(M)$.

By second countability of $M$, we may extract a countable subcover $(B_j)_{j\in\mathbb N}$ consisting of balls of the preceding type. Theorem~\ref{rademacher en variedades} implies that $\nabla|\nabla^s\mathbf{u}|_{\mathbf{g},\mathbf{h}_{\mathbf{E}}}$ exists almost everywhere on each $B_j$. Since the union is countable, it exists almost everywhere on $\operatorname{Int}(M)$ and, because $\partial M$ has measure zero, also almost everywhere on $M$.
\end{proof}

The preceding results allow us to prove Kato's inequality for sections of vector bundles.

\begin{theorem}[Kato's inequality]\label{desigualdad de kato}\glsadd{desigualdad-kato}
Let $(M,\mathbf{g})$ be a Riemannian manifold with or without boundary, and let $\mathbf{E}\to M$ be a smooth real vector bundle equipped with a bundle metric $\mathbf{h}_{\mathbf{E}}$ and a compatible connection $\nabla^{\mathbf{E}}$. For any integer $s\geq0$ and any $\mathbf{u}\in\Gamma(\mathbf{E})$, we have
$$|\nabla|\nabla^{s}\mathbf{u}|_{\mathbf{g},\mathbf{h}_{\mathbf{E}}}|_{\mathbf{g}}\leq|\nabla^{s+1}\mathbf{u}|_{\mathbf{g},\mathbf{h}_{\mathbf{E}}}$$
almost everywhere on $M$.
\end{theorem}
\begin{proof}
If $|\nabla^{s}\mathbf{u}|_{\mathbf{g},\mathbf{h}_{\mathbf{E}}}\neq0$ at a point, continuity gives a smooth chart $(U,\phi)$, with $\phi=(x^{1},\dots,x^{n})$, on which this function is positive. Thus $|\nabla^{s}\mathbf{u}|_{\mathbf{g},\mathbf{h}_{\mathbf{E}}}$ is classically differentiable on $U$, and the chain rule gives
$$\nabla|\nabla^{s}\mathbf{u}|_{\mathbf{g},\mathbf{h}_{\mathbf{E}}}=\displaystyle\frac{\nabla\langle\nabla^{s}\mathbf{u},\nabla^{s}\mathbf{u}\rangle_{\mathbf{g},\mathbf{h}_{\mathbf{E}}}}{2|\nabla^{s}\mathbf{u}|_{\mathbf{g},\mathbf{h}_{\mathbf{E}}}}.$$

To compute $\nabla\langle\nabla^{s}\mathbf{u},\nabla^{s}\mathbf{u}\rangle_{\mathbf{g},\mathbf{h}_{\mathbf{E}}}$, recall that
$$\left(\nabla\langle\nabla^{s}\mathbf{u},\nabla^{s}\mathbf{u}\rangle_{\mathbf{g},\mathbf{h}_{\mathbf{E}}}\right)_{i}=\displaystyle\frac{\partial\langle\nabla^{s}\mathbf{u},\nabla^{s}\mathbf{u}\rangle_{\mathbf{g},\mathbf{h}_{\mathbf{E}}}}{\partial x^{i}}.$$
Since the Levi--Civita connection is compatible with $\mathbf{g}$ and $\nabla^{\mathbf{E}}$ is compatible with $\mathbf{h}_{\mathbf{E}}$, the induced connection on $T^{(0,s)}(TM)\otimes \mathbf{E}$ is compatible with the product metric $\langle\cdot,\cdot\rangle_{\mathbf{g},\mathbf{h}_{\mathbf{E}}}$. Therefore,
$$\displaystyle\frac{\partial\langle\nabla^{s}\mathbf{u},\nabla^{s}\mathbf{u}\rangle_{\mathbf{g},\mathbf{h}_{\mathbf{E}}}}{\partial x^{i}}=\left\langle\nabla_{\boldsymbol{\partial}_{i}}(\nabla^{s}\mathbf{u}),\nabla^{s}\mathbf{u}\right\rangle_{\mathbf{g},\mathbf{h}_{\mathbf{E}}}+\left\langle\nabla^{s}\mathbf{u},\nabla_{\boldsymbol{\partial}_{i}}(\nabla^{s}\mathbf{u})\right\rangle_{\mathbf{g},\mathbf{h}_{\mathbf{E}}}$$
$$=2\left\langle\nabla_{\boldsymbol{\partial}_{i}}(\nabla^{s}\mathbf{u}),\nabla^{s}\mathbf{u}\right\rangle_{\mathbf{g},\mathbf{h}_{\mathbf{E}}}.$$
Thus,
$$\nabla|\nabla^{s}\mathbf{u}|_{\mathbf{g},\mathbf{h}_{\mathbf{E}}}=\displaystyle\frac{\nabla\langle\nabla^{s}\mathbf{u},\nabla^{s}\mathbf{u}\rangle_{\mathbf{g},\mathbf{h}_{\mathbf{E}}}}{2|\nabla^{s}\mathbf{u}|_{\mathbf{g},\mathbf{h}_{\mathbf{E}}}}=\displaystyle\frac{2\left\langle\nabla_{\boldsymbol{\partial}_{i}}(\nabla^{s}\mathbf{u}),\nabla^{s}\mathbf{u}\right\rangle_{\mathbf{g},\mathbf{h}_{\mathbf{E}}}\mathbf{d}x^{i}}{2|\nabla^{s}\mathbf{u}|_{\mathbf{g},\mathbf{h}_{\mathbf{E}}}}$$
$$=\displaystyle\frac{\left\langle\nabla_{\boldsymbol{\partial}_{i}}(\nabla^{s}\mathbf{u}),\nabla^{s}\mathbf{u}\right\rangle_{\mathbf{g},\mathbf{h}_{\mathbf{E}}}\mathbf{d}x^{i}}{|\nabla^{s}\mathbf{u}|_{\mathbf{g},\mathbf{h}_{\mathbf{E}}}}.$$

Consider the square root of $\mathbf{g}^{-1}$, given by Proposition \ref{raiz cuadrada simetrica}, which is positive definite and symmetric, and denote it by $(R^{ij})$.\\

\noindent Then
$$\left|\left\langle\nabla_{\boldsymbol{\partial}_{i}}(\nabla^{s}\mathbf{u}),\nabla^{s}\mathbf{u}\right\rangle_{\mathbf{g},\mathbf{h}_{\mathbf{E}}}\mathbf{d}x^{i}\right|_{\mathbf{g}}^{2}=\displaystyle\sum_{\substack{1\leq i,j\leq n}}g^{ij}\left\langle\nabla_{\boldsymbol{\partial}_{i}}(\nabla^{s}\mathbf{u}),\nabla^{s}\mathbf{u}\right\rangle_{\mathbf{g},\mathbf{h}_{\mathbf{E}}}\left\langle\nabla_{\boldsymbol{\partial}_{j}}(\nabla^{s}\mathbf{u}),\nabla^{s}\mathbf{u}\right\rangle_{\mathbf{g},\mathbf{h}_{\mathbf{E}}}$$
$$=\displaystyle\sum_{\substack{1\leq i,j\leq n}}\displaystyle\sum_{m=1}^{n}R^{im}R^{mj}\left\langle\nabla_{\boldsymbol{\partial}_{i}}(\nabla^{s}\mathbf{u}),\nabla^{s}\mathbf{u}\right\rangle_{\mathbf{g},\mathbf{h}_{\mathbf{E}}}\left\langle\nabla_{\boldsymbol{\partial}_{j}}(\nabla^{s}\mathbf{u}),\nabla^{s}\mathbf{u}\right\rangle_{\mathbf{g},\mathbf{h}_{\mathbf{E}}}$$
$$=\displaystyle\sum_{m=1}^{n}\displaystyle\sum_{1\leq i,j\leq n}R^{im}R^{jm}\left\langle\nabla_{\boldsymbol{\partial}_{i}}(\nabla^{s}\mathbf{u}),\nabla^{s}\mathbf{u}\right\rangle_{\mathbf{g},\mathbf{h}_{\mathbf{E}}}\left\langle\nabla_{\boldsymbol{\partial}_{j}}(\nabla^{s}\mathbf{u}),\nabla^{s}\mathbf{u}\right\rangle_{\mathbf{g},\mathbf{h}_{\mathbf{E}}}$$
$$=\displaystyle\sum_{m=1}^{n}\left(\displaystyle\sum_{i=1}^{n}R^{im}\left\langle\nabla_{\boldsymbol{\partial}_{i}}(\nabla^{s}\mathbf{u}),\nabla^{s}\mathbf{u}\right\rangle_{\mathbf{g},\mathbf{h}_{\mathbf{E}}}\right)^{2}$$
$$=\displaystyle\sum_{m=1}^{n}\left(\left\langle\displaystyle\sum_{i=1}^{n}R^{im}\nabla_{\boldsymbol{\partial}_{i}}(\nabla^{s}\mathbf{u}),\nabla^{s}\mathbf{u}\right\rangle_{\mathbf{g},\mathbf{h}_{\mathbf{E}}}\right)^{2}$$
$$=\displaystyle\sum_{m=1}^{n}\left|\left\langle\displaystyle\sum_{i=1}^{n}R^{im}\nabla_{\boldsymbol{\partial}_{i}}(\nabla^{s}\mathbf{u}),\nabla^{s}\mathbf{u}\right\rangle_{\mathbf{g},\mathbf{h}_{\mathbf{E}}}\right|^{2}$$
$$\underset{\text{Cauchy--Schwarz}}{\leq}\displaystyle\sum_{m=1}^{n}\left|\displaystyle\sum_{i=1}^{n}R^{im}\nabla_{\boldsymbol{\partial}_{i}}(\nabla^{s}\mathbf{u})\right|_{\mathbf{g},\mathbf{h}_{\mathbf{E}}}^{2}|\nabla^{s}\mathbf{u}|_{\mathbf{g},\mathbf{h}_{\mathbf{E}}}^{2}$$
$$=|\nabla^{s}\mathbf{u}|_{\mathbf{g},\mathbf{h}_{\mathbf{E}}}^{2}\displaystyle\sum_{m=1}^{n}\left\langle\displaystyle\sum_{i=1}^{n}R^{im}\nabla_{\boldsymbol{\partial}_{i}}(\nabla^{s}\mathbf{u}),\displaystyle\sum_{i=1}^{n}R^{im}\nabla_{\boldsymbol{\partial}_{i}}(\nabla^{s}\mathbf{u})\right\rangle_{\mathbf{g},\mathbf{h}_{\mathbf{E}}}$$
$$=|\nabla^{s}\mathbf{u}|_{\mathbf{g},\mathbf{h}_{\mathbf{E}}}^{2}\displaystyle\sum_{m=1}^{n}\left\langle\displaystyle\sum_{i=1}^{n}R^{im}\nabla_{\boldsymbol{\partial}_{i}}(\nabla^{s}\mathbf{u}),\displaystyle\sum_{j=1}^{n}R^{mj}\nabla_{\boldsymbol{\partial}_{j}}(\nabla^{s}\mathbf{u})\right\rangle_{\mathbf{g},\mathbf{h}_{\mathbf{E}}}$$
$$=|\nabla^{s}\mathbf{u}|_{\mathbf{g},\mathbf{h}_{\mathbf{E}}}^{2}\displaystyle\sum_{m=1}^{n}\left(\displaystyle\sum_{i=1}^{n}R^{im}\displaystyle\sum_{j=1}^{n}R^{mj}\left\langle\nabla_{\boldsymbol{\partial}_{i}}(\nabla^{s}\mathbf{u}),\nabla_{\boldsymbol{\partial}_{j}}(\nabla^{s}\mathbf{u})\right\rangle_{\mathbf{g},\mathbf{h}_{\mathbf{E}}}\right)$$
$$=|\nabla^{s}\mathbf{u}|_{\mathbf{g},\mathbf{h}_{\mathbf{E}}}^{2}\displaystyle\sum_{m=1}^{n}\left(\displaystyle\sum_{1\leq i,j\leq n}R^{im}R^{mj}\left\langle\nabla_{\boldsymbol{\partial}_{i}}(\nabla^{s}\mathbf{u}),\nabla_{\boldsymbol{\partial}_{j}}(\nabla^{s}\mathbf{u})\right\rangle_{\mathbf{g},\mathbf{h}_{\mathbf{E}}}\right)$$
$$=|\nabla^{s}\mathbf{u}|_{\mathbf{g},\mathbf{h}_{\mathbf{E}}}^{2}\displaystyle\sum_{1\leq i,j\leq n}\displaystyle\sum_{m=1}^{n}R^{im}R^{mj}\left\langle\nabla_{\boldsymbol{\partial}_{i}}(\nabla^{s}\mathbf{u}),\nabla_{\boldsymbol{\partial}_{j}}(\nabla^{s}\mathbf{u})\right\rangle_{\mathbf{g},\mathbf{h}_{\mathbf{E}}}$$
$$=|\nabla^{s}\mathbf{u}|_{\mathbf{g},\mathbf{h}_{\mathbf{E}}}^{2}\displaystyle\sum_{1\leq i,j\leq n}g^{ij}\left\langle\nabla_{\boldsymbol{\partial}_{i}}(\nabla^{s}\mathbf{u}),\nabla_{\boldsymbol{\partial}_{j}}(\nabla^{s}\mathbf{u})\right\rangle_{\mathbf{g},\mathbf{h}_{\mathbf{E}}}$$
$$=|\nabla^{s}\mathbf{u}|_{\mathbf{g},\mathbf{h}_{\mathbf{E}}}^{2}|\nabla(\nabla^{s}\mathbf{u})|_{\mathbf{g},\mathbf{h}_{\mathbf{E}}}^{2}=|\nabla^{s}\mathbf{u}|_{\mathbf{g},\mathbf{h}_{\mathbf{E}}}^{2}|\nabla^{s+1}\mathbf{u}|_{\mathbf{g},\mathbf{h}_{\mathbf{E}}}^{2},$$
where we used bilinearity of the induced bundle metric and the Cauchy--Schwarz inequality from Theorem~\ref{teo:cauchy-schwarz-hilbert} in each fiber of $T^{(0,s)}(TM)\otimes \mathbf{E}$. Thus,
$$|\nabla|\nabla^{s}\mathbf{u}|_{\mathbf{g},\mathbf{h}_{\mathbf{E}}}|_{\mathbf{g}}=\displaystyle\frac{\left|\left\langle\nabla_{\boldsymbol{\partial}_{i}}(\nabla^{s}\mathbf{u}),\nabla^{s}\mathbf{u}\right\rangle_{\mathbf{g},\mathbf{h}_{\mathbf{E}}}\mathbf{d}x^{i}\right|_{\mathbf{g}}}{|\nabla^{s}\mathbf{u}|_{\mathbf{g},\mathbf{h}_{\mathbf{E}}}}\leq\displaystyle\frac{|\nabla^{s}\mathbf{u}|_{\mathbf{g},\mathbf{h}_{\mathbf{E}}}|\nabla^{s+1}\mathbf{u}|_{\mathbf{g},\mathbf{h}_{\mathbf{E}}}}{|\nabla^{s}\mathbf{u}|_{\mathbf{g},\mathbf{h}_{\mathbf{E}}}}=|\nabla^{s+1}\mathbf{u}|_{\mathbf{g},\mathbf{h}_{\mathbf{E}}}.$$

Moreover, by Lemma \ref{lema para kato}, $\nabla|\nabla^{s}\mathbf{u}|_{\mathbf{g},\mathbf{h}_{\mathbf{E}}}$ exists almost everywhere on $M$. If $|\nabla^{s}\mathbf{u}|_{\mathbf{g},\mathbf{h}_{\mathbf{E}}}=0$ at a point of $\operatorname{Int}(M)$ where $\nabla|\nabla^{s}\mathbf{u}|_{\mathbf{g},\mathbf{h}_{\mathbf{E}}}$ exists, the nonnegative function $|\nabla^{s}\mathbf{u}|_{\mathbf{g},\mathbf{h}_{\mathbf{E}}}$ has a global minimum at that point. Its differential therefore vanishes, and consequently $\nabla|\nabla^{s}\mathbf{u}|_{\mathbf{g},\mathbf{h}_{\mathbf{E}}}=0$. Hence
$$|\nabla|\nabla^{s}\mathbf{u}|_{\mathbf{g},\mathbf{h}_{\mathbf{E}}}|_{\mathbf{g}}=0\leq|\nabla^{s+1}\mathbf{u}|_{\mathbf{g},\mathbf{h}_{\mathbf{E}}}.$$
We thus conclude that
$$|\nabla|\nabla^{s}\mathbf{u}|_{\mathbf{g},\mathbf{h}_{\mathbf{E}}}|_{\mathbf{g}}\leq|\nabla^{s+1}\mathbf{u}|_{\mathbf{g},\mathbf{h}_{\mathbf{E}}}$$
almost everywhere on $M$.
\end{proof}
\section{Sobolev embeddings on compact Riemannian manifolds}

Localization in coordinate balls and half-balls allows the Euclidean embeddings to be transferred to a compact Riemannian manifold. We shall use Kato's inequality to pass from first order to higher orders.

\begin{proposition}[Iteration of first-order embeddings]
\label{lema 2.11 aubin}
Let $(M,\mathbf g)$ be a Riemannian manifold with or without boundary. Let $0\leq m<k$ be integers, set $d:=k-m$, and consider exponents $p_0,\dots,p_d\in[1,\infty)$. If the inclusions
\[
W^{1,p_i}(M)\hookrightarrow L^{p_{i+1}}(M)
\qquad (i\in\{0,\dots,d-1\})
\]
are continuous, then so is
\[
W^{k,p_0}(M)\hookrightarrow W^{m,p_d}(M).
\]
\end{proposition}

\begin{proof}
We first show that a continuous embedding $W^{1,p}(M)\hookrightarrow L^q(M)$, with $p,q<\infty$, implies
\[
W^{s,p}(M)\hookrightarrow W^{s-1,q}(M)
\qquad\text{for every integer }s\geq1.
\]
Let $u\in C^{s,p}(M)$ and fix $0\leq j\leq s-1$. Kato's inequality, Theorem~\ref{desigualdad de kato}, gives
\[
\left|\nabla|\nabla^j u|_{\mathbf g}\right|_{\mathbf g}
\leq|\nabla^{j+1}u|_{\mathbf g}
\qquad\text{almost everywhere}.
\]
To work with the definition by completion, consider
\[
f_\varepsilon
:=\sqrt{|\nabla^j u|_{\mathbf g}^2+\varepsilon^2}-\varepsilon,
\qquad\varepsilon>0.
\]
These functions are smooth and satisfy $0\leq f_\varepsilon\leq|\nabla^j u|_{\mathbf g}$. Where $|\nabla^j u|_{\mathbf g}>0$, we have
\[
\nabla f_\varepsilon
=\frac{|\nabla^j u|_{\mathbf g}}
{\sqrt{|\nabla^j u|_{\mathbf g}^2+\varepsilon^2}}
\nabla|\nabla^j u|_{\mathbf g}.
\]
Where $\nabla^j u=0$, the differential of its squared norm is zero, so $\nabla f_\varepsilon=0$. Consequently,
\[
|\nabla f_\varepsilon|_{\mathbf g}
\leq|\nabla^{j+1}u|_{\mathbf g}.
\]
Moreover, $f_\varepsilon\to|\nabla^j u|_{\mathbf g}$, and their differentials converge almost everywhere to the differential of that norm: at interior zeros where it exists, the differential vanishes because the point is a minimum. Dominated convergence, applied to the $p$th powers of the preceding bounds, gives convergence of the functions in $L^p(M)$ and of their differentials in $L^p(M,T^*M)$. Thus $(f_\varepsilon)$ is Cauchy in the norm defining $W^{1,p}(M)$ and represents the function $|\nabla^j u|_{\mathbf g}$ in its completion. We also obtain
\[
\bigl\||\nabla^j u|_{\mathbf g}\bigr\|_{W^{1,p}(M)}
\leq
\|\nabla^j u\|_{L^p(M,T^{(0,j)}(TM))}
+
\|\nabla^{j+1}u\|_{L^p(M,T^{(0,j+1)}(TM))}.
\]

Apply the assumed first-order embedding to this function. There is $C_{p,q}>0$ such that
\[
\|\nabla^j u\|_{L^q(M,T^{(0,j)}(TM))}
\leq C_{p,q}\left(
\|\nabla^j u\|_{L^p(M,T^{(0,j)}(TM))}
+
\|\nabla^{j+1}u\|_{L^p(M,T^{(0,j+1)}(TM))}
\right).
\]
Summing for $j\in\{0,\dots,s-1\}$ and using Proposition~\ref{normas equivalentes de Wmp} gives
\[
\|u\|_{W^{s-1,q}(M)}\leq C_s\|u\|_{W^{s,p}(M)}.
\]
Extension to the completions, Corollary~\ref{cor:extension-operadores-completaciones}, provides the desired embedding. This extension agrees with the canonical inclusion: if a smooth sequence $(u_\nu)_{\nu\in\mathbb N}$ converges to $u$
in $L^p(M)$ and to $v$ in $L^q(M)$, then, for every compact set
$K\subseteq M$ of positive volume, Hölder's inequality gives
\[
\begin{split}
 \|u_\nu-u\|_{L^1(K)}
 &\leq\lambda_{\mathbf g}(K)^{1-1/p}\|u_\nu-u\|_{L^p(M)},\\
 \|u_\nu-v\|_{L^1(K)}
 &\leq\lambda_{\mathbf g}(K)^{1-1/q}\|u_\nu-v\|_{L^q(M)}.
\end{split}
\]
The volume of $K$ is finite. Both right-hand sides tend to zero,
so uniqueness of the limit in $L^1(K)$ implies that $u=v$ almost
everywhere on $K$. For compact sets of volume zero this last
assertion is immediate. A countable cover of $M$ by compact sets
then gives $u=v$ almost everywhere on $M$.

For each $i\in\{0,\dots,d-1\}$, apply the preceding argument with $s=k-i$, $p=p_i$, and $q=p_{i+1}$. We obtain
\[
\|u\|_{W^{k-i-1,p_{i+1}}(M)}
\leq C_i\|u\|_{W^{k-i,p_i}(M)}.
\]
The composition of these $d$ embeddings gives
\[
\|u\|_{W^{m,p_d}(M)}
\leq\left(\prod_{i=0}^{d-1}C_i\right)
\|u\|_{W^{k,p_0}(M)},
\]
as required.
\end{proof}

\begin{theorem}[Sobolev embedding with finite orders and exponents]
\label{sobolev riemanniano}
Let $(M,\mathbf g)$ be a compact Riemannian manifold of dimension $n$, with or without boundary. Let $0\leq m<k$ be integers and $1\leq p,q<\infty$ satisfy
\[
k-\frac np\geq m-\frac nq.
\]
Then the inclusion $W^{k,p}(M)\hookrightarrow W^{m,q}(M)$ is continuous.
\end{theorem}

\begin{proof}
If $n=0$, compactness implies that $M$ is a finite set, and the assertion follows from equivalence of norms in finite dimensions. Assume $n\geq1$. We begin with first order. Fix $1\leq p,q<\infty$ with $1-\displaystyle\frac np\geq-\displaystyle\frac nq$. Choose a finite cover $(U_i,\phi_i)_{i=1}^N$ consisting of interior regular coordinate balls and regular coordinate half-balls centered on the boundary. Let $(\psi_i)_{i=1}^N$ be a subordinate smooth partition of unity with $\operatorname{supp}\psi_i\Subset U_i$, and set
\[
\Omega_i:=\phi_i\bigl(U_i\cap\operatorname{Int}(M)\bigr),
\qquad
\boldsymbol\delta_i:=\phi_i^*\overline{\mathbf g},
\]
where $\overline{\mathbf g}$ is the Euclidean metric.

Let $u\in C^\infty(M)$ and define $v_i:=(\psi_i u)\circ\phi_i^{-1}$ on $\Omega_i$. If $U_i$ is interior, $v_i\in C_c^\infty(\Omega_i)$, and Theorem~\ref{sobolev generalizado}, with total orders one and zero, gives
\[
\|v_i\|_{L^q(\Omega_i)}
\leq A_i\|v_i\|_{W^{1,p}(\Omega_i)}.
\]
If $U_i$ is a boundary chart, $\Omega_i$ is a half-ball. The support of $v_i$ is separated from the spherical portion, so its extension by zero within $\mathbb R^n_+$ is smooth up to $x_n=0$ and has support in the fixed compact set
\[
K_i:=\phi_i(\operatorname{supp}\psi_i)
\subseteq\overline{\mathbb R^n_+}.
\]
Theorem~\ref{teo:encajes-localizados-semiespacio}, with total orders one and zero, gives the same estimate upon restriction to $\Omega_i$.

Change of variables and the definition of the metric $\boldsymbol\delta_i$ give
\[
\int_{\Omega_i}|v_i|^q\,d\lambda_n
=\int_{U_i}|\psi_i u|^q\,d\lambda_{\boldsymbol\delta_i},
\qquad
\int_{\Omega_i}|dv_i|_{\overline{\mathbf g}}^p\,d\lambda_n
=\int_{U_i}|d(\psi_i u)|_{\boldsymbol\delta_i}^p\,
d\lambda_{\boldsymbol\delta_i}.
\]
Equivalence of norms on $\mathbb R^n$, Proposition~\ref{normas en R^{n}}, compares $|dv_i|_{\overline{\mathbf g}}$ with $\bigl(\displaystyle\sum_{a=1}^n|D_a v_i|^p\bigr)^{1/p}$. On the other hand, Remark~\ref{obs:comparacion-metricas-compactos-coordenados} allows comparison of the metrics $\boldsymbol\delta_i$ and $\mathbf g$ and their measures on $\overline{U_i}$. Thus there exist $a_i,b_i>0$, independent of $u$, such that
\[
\|v_i\|_{W^{1,p}(\Omega_i)}
\leq a_i\|\psi_i u\|_{W^{1,p}(U_i),\mathbf g},
\qquad
\|\psi_i u\|_{L^q(U_i),\mathbf g}
\leq b_i\|v_i\|_{L^q(\Omega_i)}.
\]
The Euclidean estimate then implies
\[
\|\psi_i u\|_{L^q(M)}
\leq a_i b_i A_i\|\psi_i u\|_{W^{1,p}(M)},
\]
since $\psi_i u$ and its derivatives vanish outside $U_i$.

To bound the localization, use $d(\psi_i u)=\psi_i\,du+u\,d\psi_i$ and $0\leq\psi_i\leq1$. Since $M$ is compact, $|d\psi_i|_{\mathbf g}$ is bounded and
\[
\begin{aligned}
\|\psi_i u\|_{W^{1,p}(M)}
&\leq\|\psi_i u\|_{L^p(M)}
+\|d(\psi_i u)\|_{L^p(M,T^*M)}\\
&\leq
\bigl(1+\||d\psi_i|_{\mathbf g}\|_{L^\infty(M)}\bigr)
\|u\|_{L^p(M)}
+\|du\|_{L^p(M,T^*M)}\\
&\leq C_i\|u\|_{W^{1,p}(M)}.
\end{aligned}
\]
Summing over the finite cover and using $u=\displaystyle\sum_{i=1}^N\psi_i u$ gives
\[
\|u\|_{L^q(M)}
\leq\sum_{i=1}^N\|\psi_i u\|_{L^q(M)}
\leq C\|u\|_{W^{1,p}(M)}.
\]
The density defining $W^{1,p}(M)$ and Theorem~\ref{teo:extension-operadores-subespacio-denso} extend this estimate to all of $W^{1,p}(M)$. The limit in $L^q(M)$ is identified with the limit in $L^p(M)$ by uniqueness of the limit in $L^1(M)$.

Now consider the orders $m<k$ and the exponents in the statement. Set $d:=k-m$ and define
\[
\frac1{p_i}
:=\frac1p-\frac{i}{d}\left(\frac1p-\frac1q\right),
\qquad i\in\{0,\dots,d\}.
\]
Each $p_i$ is finite and belongs to $[1,\infty)$, since its reciprocal is a convex combination of $1/p$ and $1/q$. Moreover, $p_0=p$, $p_d=q$, and
\[
1-\frac n{p_i}+\frac n{p_{i+1}}
=1-\frac nd\left(\frac1p-\frac1q\right)\geq0.
\]
The first-order part gives $W^{1,p_i}(M)\hookrightarrow L^{p_{i+1}}(M)$ continuously for each $i\in\{0,\dots,d-1\}$. Proposition~\ref{lema 2.11 aubin} then yields $W^{k,p}(M)\hookrightarrow W^{m,q}(M)$ continuously.
\end{proof}

\begin{corollary}[Sobolev embedding and critical exponent]
\label{sobolev riemanniano 2}
Let $(M,\mathbf g)$ be a compact Riemannian manifold of dimension $n$, with or without boundary. Let $j\geq0$ and $d\geq1$ be integers, and let $1\leq p<n/d$. Then
\[
W^{j+d,p}(M)\hookrightarrow W^{j,q}(M)
\qquad\text{for every }1\leq q\leq
p_d^*:=\frac{np}{n-dp},
\]
and these inclusions are continuous.
\end{corollary}

\begin{proof}
The condition $q\leq p_d^*$ is equivalent to $j+d-\displaystyle\frac np\geq j-\displaystyle\frac nq$. Apply Theorem~\ref{sobolev riemanniano} with $k=j+d$ and $m=j$.
\end{proof}

\section{The Rellich--Kondrashov theorem on compact Riemannian manifolds}
The Rellich--Kondrashov theorem provides the compact counterpart of the preceding embeddings. The Euclidean version in Theorem~\ref{rellich kondrashov generalizado} can be transferred to a compact Riemannian manifold $(M,\mathbf{g})$ by means of a finite cover by charts and a partition of unity; in boundary charts, we use the localized half-space version. The following result specifies the orders and exponents for which compactness holds, including the loss of one derivative with any finite exponent. In the range where a finite critical exponent exists, we shall recover the subcritical embeddings from Corollary~\ref{sobolev riemanniano 2} as a special case.

\begin{theorem}[Rellich--Kondrashov]\label{teo:sobolev-en-variedades-de-rellich-kondrashov} Let $(M,\mathbf{g})$ be a compact Riemannian manifold, with or without boundary, of dimension $n$. For any integers $j\geq0$ and $m\geq1$, and any $1\leq p,q<\infty$ satisfying $j+m-\frac np>j-\frac nq$, the inclusion $W^{j+m,p}(M)\subseteq W^{j,q}(M)$ is compact.
\end{theorem}

The proof uses the Leibniz rule for the tensor connection, Proposition~\ref{propiedades de la derivada covariante tensorial}, iterated while preserving the order of the covariant indices.

\begin{proof}
Since $M$ is compact, there is a finite open cover $\{(U_{i},\phi_{i})\mid i\in\{1,\dots,N\}\}$ consisting of regular coordinate balls contained in $\operatorname{Int}(M)$ and, when $\partial M\neq\varnothing$, regular coordinate half-balls centered on the boundary. Let $(\psi_{i})_{i=1}^{N}$ be a partition of unity subordinate to this open cover. In both cases, $\phi_{i}$ extends to $\overline{U_{i}}$, which is compact. Consider again $\boldsymbol{\delta}_{i}=\phi_{i}^{*}\overline{\mathbf{g}}$.\\

Let $(u_{k})_{k\in\mathbb{N}}\subseteq W^{j+m,p}(M)\cap C^{\infty}(M)$ be a bounded sequence.\\

By Remark~\ref{obs:comparacion-metricas-compactos-coordenados}, the local estimates in Theorem~\ref{W no depende de metrica} apply on each $\overline{U_i}$, even when a half-ball has a corner. Thus, for each $i\in\{1,\dots,N\}$, there exist $C_{i},c_{i}>0$ such that $c_{i}\|\psi_{i}u_{k}\|_{W^{j+m,p}(\overline{U_{i}}),\boldsymbol{\delta}_{i}}\leq\|\psi_{i}u_{k}\|_{W^{j+m,p}(\overline{U_{i}}),\mathbf{g}}\leq C_{i}\|\psi_{i}u_{k}\|_{W^{j+m,p}(\overline{U_{i}}),\boldsymbol{\delta}_{i}}$, where the constants are independent of the function. Set $v_{i,k}:=(\psi_i u_k)\circ\phi_i^{-1}$. Change of variables gives
\[
\|\psi_i u_k\|_{W^{j+m,p}(U_i),\mathbf g}
\geq c_i\left(
\sum_{s=0}^{j+m}\int_{\phi_i(U_i)}
|\nabla^s v_{i,k}|_{\overline{\mathbf g}}^p\,d\lambda_n
\right)^{\frac1p}.
\]
For each $s\in\{0,\dots,j+m\}$, Proposition~\ref{normas en R^{n}} provides $M_s,m_s>0$ such that $m_s\|\cdot\|_p\leq\|\cdot\|_2\leq M_s\|\cdot\|_p$ on $\mathbb R^{n^s}$, where $n^s$ counts the ordered blocks of $s$ indices. Each multi-index $\alpha$ with $|\alpha|=s$ occurs $s!/\alpha!$ times among these blocks, with $\alpha!:=\alpha_1!\cdots\alpha_n!$. Therefore,
\[
\begin{aligned}
|\nabla^s v_{i,k}|_{\overline{\mathbf g}}^p
&=\left(\sum_{|\alpha|=s}\frac{s!}{\alpha!}
|D^\alpha v_{i,k}|^2\right)^{\frac p2}\\
&\geq m_s^p\sum_{|\alpha|=s}|D^\alpha v_{i,k}|^p.
\end{aligned}
\]
Integrating and summing over the orders gives
\[
\begin{aligned}
\sum_{s=0}^{j+m}\int_{\phi_i(U_i)}
|\nabla^s v_{i,k}|_{\overline{\mathbf g}}^p\,d\lambda_n
&\geq
\sum_{s=0}^{j+m}m_s^p
\sum_{|\alpha|=s}\int_{\phi_i(U_i)}
|D^\alpha v_{i,k}|^p\,d\lambda_n\\
&\geq
(m')^p\sum_{|\alpha|\leq j+m}
\int_{\phi_i(U_i)}|D^\alpha v_{i,k}|^p\,d\lambda_n,
\end{aligned}
\]
where $m':=\min_{0\leq s\leq j+m}m_s>0$. Taking the $p$th root and using the definition of the Euclidean Sobolev norm, we conclude that
\[
\|\psi_i u_k\|_{W^{j+m,p}(U_i),\mathbf g}
\geq c_i m'\|v_{i,k}\|_{W^{j+m,p}(\phi_i(U_i))}.
\]

The iterated Leibniz rule expresses $\nabla^s(\psi_i u_k)$ as a sum of terms obtained by distributing the $s$ derivatives between $\psi_i$ and $u_k$. For each $r\in\{0,\dots,s\}$, there are $\binom sr$ terms with $s-r$ derivatives on $\psi_i$ and $r$ on $u_k$. These terms differ by permutations of the covariant indices, which preserve the tensor norm. Proposition~\ref{norma producto tensorial} gives
\[
|\mathbf F\otimes\mathbf G|_{\mathbf g}
=|\mathbf F|_{\mathbf g}|\mathbf G|_{\mathbf g}
=|\mathbf G\otimes\mathbf F|_{\mathbf g}.
\]
Thus, for $0\leq s\leq j+m$,
\[
\begin{aligned}
|\nabla^s(\psi_i u_k)|_{\mathbf g}
&\leq
\sum_{r=0}^s\binom sr
|\nabla^{s-r}\psi_i|_{\mathbf g}|\nabla^r u_k|_{\mathbf g}\\
&\leq
\sum_{r=0}^s\binom sr
\||\nabla^{s-r}\psi_i|_{\mathbf g}\|_{L^\infty(M)}
|\nabla^r u_k|_{\mathbf g}.
\end{aligned}
\]
Let
\[
A_i:=\max_{0\leq s\leq j+m}
\||\nabla^s\psi_i|_{\mathbf g}\|_{L^\infty(M)}.
\]
Compactness of $M$ implies $A_i<\infty$. Integrating the pointwise estimate and applying Minkowski's inequality gives
\[
\|\nabla^s(\psi_i u_k)\|_{L^p(M,T^{(0,s)}(TM))}
\leq A_i\sum_{r=0}^s\binom sr
\|\nabla^r u_k\|_{L^p(M,T^{(0,r)}(TM))}.
\]
Now sum over $s$ and use Proposition~\ref{normas equivalentes de Wmp}:
\[
\begin{aligned}
\|\psi_i u_k\|_{W^{j+m,p}(M)}
&\leq A_i\sum_{s=0}^{j+m}\sum_{r=0}^s\binom sr
\|\nabla^r u_k\|_{L^p(M,T^{(0,r)}(TM))}\\
&\leq C_i'\|u_k\|_{W^{j+m,p}(M)}.
\end{aligned}
\]
The constant $C_i'$ is independent of $k$. Together with the preceding comparison, this gives
\[
\frac{c_i m'}{C_i'}\|v_{i,k}\|_{W^{j+m,p}(\phi_i(U_i))}
\leq\|u_k\|_{W^{j+m,p}(M)}.
\]
Thus $(v_{i,k})_k$ is bounded in the corresponding Euclidean space. If $U_i$ is interior, the support of each $v_{i,k}$ is contained in the fixed compact set $\phi_i(\operatorname{supp}\psi_i)\Subset\phi_i(U_i)$. Theorem~\ref{rellich kondrashov generalizado} provides a subsequence converging in $W^{j,q}(\phi_i(U_i))$.

If $U_i$ is a boundary chart, extend $v_{i,k}$ by zero outside the half-ball, but within $\mathbb R^n_+$. Lemma~\ref{lem:extension-cero-media-bola-semiespacio} preserves its Sobolev norm, and the extensions have support contained in the fixed compact set $K_i:=\phi_i(\operatorname{supp}\psi_i)
\subseteq\overline{\mathbb R^n_+}$. Theorem~\ref{teo:encajes-localizados-semiespacio} gives a subsequence converging in $W^{j,q}(\mathbb R^n_+)$ and, upon restriction, a Cauchy subsequence in the Sobolev space on the half-ball.

Since there are finitely many charts, extract successively a subsequence for $U_1$, then a subsequence of that one for $U_2$, and so on up to $U_N$. Relabeling the last subsequence, we may assume without loss of generality that $(\psi_{i}u_{k}\circ\phi_{i}^{-1})_{k\in\mathbb{N}}$ is Cauchy in $W^{j,q}(\phi_{i}(U_{i}))$ for every $i\in\{1,\dots,N\}$. We shall show that $(\psi_{i}u_{k})_{k\in\mathbb{N}}$ is Cauchy in $W^{j,q}(U_{i})$.\\

Repeating a procedure analogous to the preceding one, if $\widetilde{C_{i}},\widetilde{c_{i}}>0$ satisfy $$\widetilde{c}_{i}\|\psi_{i}u_{k}\|_{W^{j,q}(\overline{U_{i}}),\boldsymbol{\delta}_{i}}\leq\|\psi_{i}u_{k}\|_{W^{j,q}(\overline{U_{i}}),\mathbf{g}}\leq \widetilde{C}_{i}\|\psi_{i}u_{k}\|_{W^{j,q}(\overline{U_{i}}),\boldsymbol{\delta}_{i}}$$, where the constants are independent of the function, and $\widetilde{M}'=\displaystyle\max_{0\leq s\leq j}\bigl(\widetilde{M}_{s}(s!)^{1/q}\bigr)$, where $\widetilde{M}_{s},\widetilde{m}_{s}>0$ satisfy $\widetilde{m}_{s}\|\cdot\|_{q}\leq \|\cdot\|_{2}\leq\widetilde{M}_{s}\|\cdot\|_{q}$ on $\mathbb{R}^{n^{s}}$ with $s\in\{0,\dots, j\}$, use $s!/\alpha!\leq s!$ to compare the sum over blocks with the sum over multi-indices. This includes $j=0$, with $0!=1$. Then $$\|\psi_{i}u_{k}-\psi_{i}u_{l}\|_{W^{j,q}(U_{i}),\mathbf{g}}\leq \widetilde{C}_{i}\widetilde{M}'\|\psi_{i}u_{k}\circ\phi_{i}^{-1}-\psi_{i}u_{l}\circ\phi_{i}^{-1}\|_{W^{j,q}(\phi_{i}(U_{i})),\overline{\mathbf{g}}}.$$

Thus, the sequence $(\psi_{i}u_{k})_{k\in\mathbb{N}}$ is Cauchy in $W^{j,q}(U_{i})$. Finally, if $k,l\in\mathbb{N}$, then
$$\|u_{k}-u_{l}\|_{W^{j,q}(M),\mathbf{g}}=\left\|\displaystyle\sum_{i=1}^{N}\psi_{i}(u_{k}-u_{l})\right\|_{W^{j,q}(M),\mathbf{g}}\leq \displaystyle\sum_{i=1}^{N}\|\psi_{i}(u_{k}-u_{l})\|_{W^{j,q}(M),\mathbf{g}}$$
$$=\displaystyle\sum_{i=1}^{N}\|\psi_{i}(u_{k}-u_{l})\|_{W^{j,q}(U_{i}),\mathbf{g}}=\displaystyle\sum_{i=1}^{N}\|\psi_{i}u_{k}-\psi_{i}u_{l}\|_{W^{j,q}(U_{i}),\mathbf{g}}.$$
Hence $(u_{k})_{k\in\mathbb{N}}$ is Cauchy in $W^{j,q}(M)$. We have proved that the restricted inclusion $I_{0}\colon C^{j+m,p}(M)\longrightarrow W^{j,q}(M)$, given by $I_{0}u=u$, is compact. Moreover, $I_0$ is continuous: the compactness just proved maps its unit ball into a relatively compact, and hence bounded, subset of $W^{j,q}(M)$; linearity turns this bound into continuity. Since $C^{j+m,p}(M)$ is dense in $W^{j+m,p}(M)$ by definition and $W^{j,q}(M)$ is a Banach space, Theorem~\ref{teo:extension-operadores-subespacio-denso} provides a unique compact linear extension $I\colon W^{j+m,p}(M)\longrightarrow W^{j,q}(M)$. This extension agrees with the canonical inclusion. Therefore, the inclusion $W^{j+m,p}(M)\subseteq W^{j,q}(M)$ is compact for all the orders and exponents in the statement.
\end{proof}

\begin{corollary}[Total orders and loss of one derivative]
\label{cor:rellich-general-variedades-compactas}
Let $(M,\mathbf{g})$ be a compact Riemannian manifold of dimension $n$, with possibly empty smooth boundary. If $m>k\geq0$ are integers and $1\leq p,q<\infty$ satisfy $m-\frac np>k-\frac nq$, then $W^{m,p}(M)\hookrightarrow W^{k,q}(M)$ is compact. In particular,
\[
W^{m,p}(M)\hookrightarrow W^{m-1,p}(M)
\]
is compact for every $m\geq1$ and every $1\leq p<\infty$.
\end{corollary}
\begin{proof}
Apply the preceding theorem with $j=k$ and a gain of $m-k\geq1$ derivatives. For the final embedding, take $k=m-1$ and $q=p$; the strict inequality reduces to $1>0$. Thus this case includes $p=1$ and dimensions $n=1,2$ without introducing a critical exponent.
\end{proof}

\begin{corollary}[Subcritical formulation]
\label{cor:rellich-subcritico-variedades-compactas}
Under the geometric assumptions of the preceding corollary, if $j\geq0$, $a\geq1$, and $ap<n$, the inclusion $W^{j+a,p}(M)\hookrightarrow W^{j,q}(M)$ is compact for $1\leq q<np/(n-ap)$. At the endpoint $q=np/(n-ap)$, Corollary~\ref{sobolev riemanniano 2} gives continuity.
\end{corollary}
\begin{proof}
With $ap<n$, the inequality $j+a-\frac np>j-\frac nq$ is equivalent to $q<np/(n-ap)$, so the general theorem applies. Equality of indices corresponds to the critical exponent in the cited continuous embedding theorem.
\end{proof}

\chapter{Partial differential operators on vector bundles}
\label{cap:operadores-diferenciales-parciales-haces}
Differential operators on manifolds naturally act between sections of vector bundles. The gradient, exterior derivative, connections, Hodge operators, and geometric Laplacians are examples of a common structure: their value at a point depends only on finitely many local derivatives of the section.

The intrinsic formulation will allow us to study the order of a differential operator, composition, and commutators with smooth functions without favoring coordinates. The formal adjoint will provide the integration by parts identity needed to define weak covariant derivatives and bundle-valued distributions later. The principal symbol and ellipticity will appear once these two theories have been constructed, so that their consequences can be expressed directly in terms of Sobolev sections and geometric distributions.

\begin{semblanzaHistorica}{From local operators to the language of bundles}
For a long time, differential operators were presented through coordinate formulas. Differential geometry showed that the essential feature was not the particular list of coefficients, but the local dependence on finitely many derivatives and the way those expressions change between trivializations. The language of bundles brings operators as seemingly different as the exterior derivative, connections, and Dirac operators together in one theory. Later, the principal symbol will retain precisely the highest-order part of this local information.
\end{semblanzaHistorica}

\section{Motivation and basic properties}

Before turning to the intrinsic definition, it is helpful to recall the Euclidean case. There, operators are written in terms of partial derivatives of the components; under a change of coordinates, those expressions transform consistently and determine a global map between spaces of sections. The abstract definition will isolate precisely this local dependence on finitely many derivatives.

Consider the Laplacian
\[
\Delta\colon C^{\infty}(\mathbb{R}^{n}) \longrightarrow C^{\infty}(\mathbb{R}^{n}),
\qquad \Delta u = \displaystyle\sum_{i=1}^{n}\frac{\partial^{2} u}{\partial x_{i}^{2}}.
\]

This operator can be interpreted as
\[
\Delta\colon \Gamma(\underline{\mathbb{R}})\longrightarrow\Gamma(\underline{\mathbb{R}}),
\]
where $\underline{\mathbb{R}}$ denotes the trivial bundle of rank $1$ over $\mathbb{R}^{n}$, that is,
\[
\pi\colon \mathbb{R}^{n}\times \mathbb{R}\longrightarrow \mathbb{R}^{n},\qquad \pi(x,v)=x.
\]

Another example is the exterior derivative
\[
d\colon \Omega^{k}(\mathbb{R}^{n}) \longrightarrow \Omega^{k+1}(\mathbb{R}^{n}),
\]
where a $k$-form $\boldsymbol{\omega}\in \Omega^{k}(\mathbb{R}^{n})$ has, in the standard coordinates, a smooth representative
\[
\widetilde\omega\colon \mathbb{R}^{n} \longrightarrow \mathbb{R}^{{n\choose k}}.
\]
Thus, $d$ can be written as an operator
\[
d\colon C^{\infty}\left(\mathbb{R}^{n},\mathbb{R}^{{n\choose k}}\right)\longrightarrow C^{\infty}\left(\mathbb{R}^{n},\mathbb{R}^{{n\choose k+1}}\right).
\]

\begin{remark}\label{obs:operadores-diferenciales-en-haces-caso-haz-trivial-rango-denotado-seccion}
In the case of the trivial bundle of rank $m$ over $\mathbb{R}^n$, denoted by $\underline{\mathbb{R}}^m$,
\[
\pi\colon \mathbb{R}^n\times\mathbb{R}^m\longrightarrow\mathbb{R}^n,
\]
every smooth section $\boldsymbol{\sigma}$ has the form $\boldsymbol{\sigma}(x)=(x,f(x))$ for a unique function $f\colon \mathbb{R}^n\longrightarrow\mathbb{R}^m$. This gives a natural identification
\[
\Gamma(\underline{\mathbb{R}}^m)\cong C^\infty(\mathbb{R}^n,\mathbb{R}^m).
\]

More generally, if $M$ is a smooth manifold with or without boundary and we consider the trivial bundle of rank $m$, again denoted by $\underline{\mathbb{R}}^m$,
\[
\pi\colon M\times\mathbb{R}^m\longrightarrow M,
\]
the map sending a section $\boldsymbol{\sigma}(p)=(p,f(p))$ to its second component $f$ is linear and bijective, with inverse $f\mapsto[p\mapsto(p,f(p))]$. Therefore,
\[
\Gamma(\underline{\mathbb{R}}^m)\cong C^\infty(M,\mathbb{R}^m).
\]
\end{remark}

In this sense, the exterior derivative may be viewed as
\[
d\colon \Gamma\left(\underline{\mathbb{R}}^{{n\choose k}}\right)\longrightarrow
\Gamma\left(\underline{\mathbb{R}}^{{n\choose k+1}}\right).
\]

These examples show that many classical differential operators can naturally be described as operators between spaces of sections of vector bundles. This setting can be generalized as follows:

\begin{definition}\label{def:operadores-diferenciales-en-haces-op}\index{space of linear operators between sections}
 Let $M$ be a smooth manifold with or without boundary, and let $\pi_{1}\colon \mathbf{E}\longrightarrow M$, $\pi_{2}\colon \mathbf{F}\longrightarrow M$ be smooth vector bundles (real or complex) over $M$. We write \[\mathbf{Op}(\mathbf{E},\mathbf{F}):=\{ P\colon \Gamma(\mathbf{E})\longrightarrow \Gamma(\mathbf{F})\mid \text{$P$ is $\mathbb{K}$-linear}\},\] where $\mathbb{K}\in \{\mathbb{R},\mathbb{C}\}$.
\end{definition}
The spaces of smooth sections $\Gamma(\mathbf{E})$ and $\Gamma(\mathbf{F})$ are not only $\mathbb{K}$-vector spaces, but also $C^{\infty}(M)$-modules. Partial differential operators are elements of $\operatorname{\mathbf{Op}}$ that interact in a special way with the $C^{\infty}(M)$-module structure of the spaces of sections.

An element $P\in \mathbf{Op}(\mathbf{E},\mathbf{F})$ is $C^{\infty}(M)$-linear if $P(f\mathbf{u})=f(P \mathbf{u})$ for any $\mathbf{u}\in \Gamma(\mathbf{E})$. We define the commutator of $P\in \mathbf{Op}(\mathbf{E},\mathbf{F})$ and $f\in C^{\infty}(M)$ to be the element $[P,f]\in \mathbf{Op}(\mathbf{E},\mathbf{F})$ such that $[P,f](\mathbf{u})=P(f\mathbf{u})-f(P\mathbf{u})$ for every $\mathbf{u}\in \Gamma(\mathbf{E})$.

The operator $P$ is $C^{\infty}(M)$--linear if and only if $[P,f]=0$ for any $f\in C^{\infty}(M)$.

Each function $f\in C^{\infty}(M)$ defines a map
\[
\text{ad}(f)\colon \mathbf{Op}(\mathbf{E},\mathbf{F})\longrightarrow \mathbf{Op}(\mathbf{E},\mathbf{F}),
\qquad \text{ad}(f)(P)=[P,f].
\]

We define a partial differential operator recursively as follows:

\begin{definition}\label{def:operadores-diferenciales-en-haces-pdo}\index{partial differential operator}
 Let $M$ be a smooth manifold with or without boundary and let $\mathbf{E},\mathbf{F}\longrightarrow M$ be smooth vector bundles over $M$.
 \begin{enumerate}[label=(\alph*)]
 \item $\mathbf{PDO}^{(0)}(\mathbf{E},\mathbf{F}):=\{P\in \mathbf{Op}(\mathbf{E},\mathbf{F})\mid [P,f]=0,\forall f\in C^\infty(M)\}$.
 \item If $m>0$,
 \[
 \mathbf{PDO}^{(m)}(\mathbf{E},\mathbf{F}):=\{P\in \mathbf{Op}(\mathbf{E},\mathbf{F})\mid [P,f]\in \mathbf{PDO}^{(m-1)}(\mathbf{E},\mathbf{F}),\forall f\in C^{\infty}(M)\}.
 \]
 \item $\mathbf{PDO}(\mathbf{E},\mathbf{F}):=\displaystyle\bigcup_{m=0}^{\infty}\mathbf{PDO}^{(m)}(\mathbf{E},\mathbf{F})$.
 \end{enumerate}
 In general, if $P\in \mathbf{PDO}^{(m)}(\mathbf{E},\mathbf{F})$, we say that $P$ is a partial differential operator of order less than or equal to $m$.
 \end{definition}

In terms of the notation $\mathrm{ad}$, this definition can be rewritten more explicitly:

\begin{proposition}\label{equivalencia pdo con ad}
 Let $M$ be a smooth manifold with or without boundary and let $\mathbf{E},\mathbf{F}\longrightarrow M$ be smooth vector bundles over $M$.
 \begin{enumerate}[label=(\alph*)]
 \item $\mathbf{PDO}^{(0)}(\mathbf{E},\mathbf{F})=\displaystyle\bigcap_{f\in C^\infty(M)}\ker(\mathrm{ad}(f))$.
 \item More generally, for $m\geq 1$,
 \[
 \mathbf{PDO}^{(m)}(\mathbf{E},\mathbf{F})=
 \displaystyle\bigcap_{f_0,\dots,f_m\in C^\infty(M)}
 \ker\big(\mathrm{ad}(f_0)\cdots \mathrm{ad}(f_m)\big).
 \]
 In particular, $P\in \mathbf{PDO}^{(m)}(\mathbf{E},\mathbf{F})$ if and only if
 \[
 \mathrm{ad}(f_0)\cdots \mathrm{ad}(f_m)(P)=0,
 \qquad \forall f_0,\dots,f_m\in C^\infty(M).
 \]
 \end{enumerate}
\end{proposition}
\begin{proof}
For the base case \(m=0\) we use the definition; note that \(P\in \mathbf{PDO}^{(0)}(\mathbf{E},\mathbf{F})\) if and only if \([P,f]=0\) for every \(f\in C^{\infty}(M)\), that is, \(\operatorname{ad}(f)(P)=0\) for every \(f\in C^{\infty}(M)\), which is equivalent to
\[
P\in \displaystyle\bigcap_{f\in C^{\infty}(M)}\ker(\operatorname{ad}(f)).
\]

Now suppose that the characterization holds up to \(m-1\), that is,
\[
\mathbf{PDO}^{(m-1)}(\mathbf{E},\mathbf{F})
=\displaystyle\bigcap_{f_0,\dots,f_{m-1}\in C^\infty(M)}
\ker\big(\operatorname{ad}(f_0)\cdots \operatorname{ad}(f_{m-1})\big).
\]
We wish to show that the same description holds for \(m\).

By definition,
\[
P\in \mathbf{PDO}^{(m)}(\mathbf{E},\mathbf{F})
\quad\Longleftrightarrow\quad
[P,f_m]\in \mathbf{PDO}^{(m-1)}(\mathbf{E},\mathbf{F})
\text{ for every } f_m\in C^\infty(M).
\]
Applying the induction hypothesis, this condition is equivalent to
\[
\operatorname{ad}(f_0)\cdots \operatorname{ad}(f_{m-1})([P,f_m])=0,
\qquad
\forall f_0,\dots,f_{m-1},f_m\in C^\infty(M).
\]
But this holds if and only if
\[
\operatorname{ad}(f_0)\cdots
\operatorname{ad}(f_{m-1})\operatorname{ad}(f_m)(P)=0,
\qquad
\forall f_0,\dots,f_m\in C^\infty(M),
\]
that is, if and only if
\[
P\in
\displaystyle\bigcap_{f_0,\dots,f_m\in C^\infty(M)}
\ker\big(\operatorname{ad}(f_0)\cdots \operatorname{ad}(f_m)\big).
\]
Therefore, the equivalence holds for \(m\), completing the induction.
\end{proof}

\begin{remark}\label{obs:operadores-diferenciales-en-haces-todas}
In the literature, one commonly encounters the abbreviation
\[
\mathbf{PDO}^{(m)}(\mathbf{E},\mathbf{F}) = \ker \mathrm{ad}^{m+1},
\]
which is to be understood in the preceding sense: not as a power of a single operator, but as the intersection of the kernels of \emph{all} possible compositions of $m+1$ operators of the form $\mathrm{ad}(f)$ with $f\in C^\infty(M)$.
\end{remark}

\begin{proposition}\label{prop:propiedades-conmutador}
Let $M$ be a smooth manifold with or without boundary and let $\mathbf{E},\mathbf{F}\longrightarrow M$ be smooth vector bundles. Let $P,Q\in \operatorname{\mathbf{Op}}(\mathbf{E},\mathbf{F})$, $f,g,h\in C^{\infty}(M)$, and $\alpha,\beta\in \mathbb{K}$. The commutator
\[
[P,f](\mathbf{u}):=P(f\mathbf{u})-fP(\mathbf{u}),\qquad \mathbf{u}\in \Gamma(\mathbf{E}),
\]
satisfies the following properties:
\begin{enumerate}[label=(\alph*)]
 \item $\mathbb{K}$-linearity in the operator:
 \[
 [\alpha P+\beta Q,f]
 =\alpha\, [P,f]+\beta\, [Q,f].
 \]

 \item $C^{\infty}(M)$-linearity in the operator:
 \[
 [hP,f]=h\,[P,f].
 \]

 \item $\mathbb{K}$-linearity in the function:
 \[
 [P,\alpha f+\beta g]
 =\alpha\, [P,f]+\beta\, [P,g].
 \]

 \item Leibniz rule:
 for every section $\mathbf{u}\in\Gamma(\mathbf{E})$,
 \[
 [P,fg](\mathbf{u})
 =[P,f](g\mathbf{u})+f\, [P,g](\mathbf{u}).
 \]
 If $P\in\mathbf{PDO}^{(1)}(\mathbf{E},\mathbf{F})$, this identity reduces to
 \[
 [P,fg](\mathbf{u})
 =g\, [P,f](\mathbf{u})+f\, [P,g](\mathbf{u}).
 \]

 \item If $c\in\mathbb{K}$ is viewed as a constant function, then
 \[
 [P,c]=0.
 \]

 \item The operator $P$ is $C^{\infty}(M)$-linear if and only if
 \[
 [P,f]=0\qquad \forall f\in C^{\infty}(M).
 \]
\end{enumerate}
\end{proposition}

\begin{proof}
Let $\mathbf{u}\in \Gamma(\mathbf{E})$.

(a) We have
\[
[\alpha P+\beta Q,f](\mathbf{u})
=(\alpha P+\beta Q)(f\mathbf{u})-f(\alpha P+\beta Q)(\mathbf{u})
=\alpha(P(f\mathbf{u})-fP(\mathbf{u}))+\beta(Q(f\mathbf{u})-fQ(\mathbf{u})),
\]
which is $\alpha [P,f](\mathbf{u})+\beta [Q,f](\mathbf{u})$.

(b) By definition,
\[
[hP,f](\mathbf{u})
=(hP)(f\mathbf{u})-f(hP)(\mathbf{u})
=hP(f\mathbf{u})-fhP(\mathbf{u})
=h\bigl(P(f\mathbf{u})-fP(\mathbf{u})\bigr)
=h\,[P,f](\mathbf{u}).
\]

(c) For $\alpha f+\beta g$,
\[
[P,\alpha f+\beta g](\mathbf{u})
=P((\alpha f+\beta g)\mathbf{u})-(\alpha f+\beta g)P(\mathbf{u})
=\alpha(P(f\mathbf{u})-fP(\mathbf{u}))+\beta(P(g\mathbf{u})-gP(\mathbf{u})).
\]

(d) For $fg$,
\[
[P,fg](\mathbf{u})=P(fg\,\mathbf{u})-fgP(\mathbf{u}).
\]
We add and subtract $fP(g\mathbf{u})$ to obtain
\[
[P,fg](\mathbf{u})
=\bigl(P(fg\mathbf{u})-fP(g\mathbf{u})\bigr)
+f\bigl(P(g\mathbf{u})-gP(\mathbf{u})\bigr)
=[P,f](g\mathbf{u})+f[P,g](\mathbf{u}).
\]
If $P\in\mathbf{PDO}^{(1)}(\mathbf{E},\mathbf{F})$, then $[P,f]\in\mathbf{PDO}^{(0)}(\mathbf{E},\mathbf{F})$ by definition. In particular, $[[P,f],g]=0$, that is,
\[
[P,f](g\mathbf{u})=g[P,f](\mathbf{u}).
\]
Substituting this equality gives the formula for the first-order case.

(e) If $c$ is constant, then
\[
[P,c](\mathbf{u})=P(c\mathbf{u})-cP(\mathbf{u})=cP(\mathbf{u})-cP(\mathbf{u})=0.
\]

(f)
If $P$ is $C^{\infty}(M)$-linear, then $P(f\mathbf{u})=fP(\mathbf{u})$, and hence $[P,f]=0$. Conversely, if $[P,f]=0$ for every $f$, then
\[
P(f\mathbf{u})-fP(\mathbf{u})=0,
\]
for every section $\mathbf{u}$, implying that $P$ is $C^{\infty}(M)$-linear.
\end{proof}

\begin{proposition}\label{prop:propiedades-ad}
Let $M$ be a smooth manifold with or without boundary and let $\mathbf{E},\mathbf{F}\longrightarrow M$ be smooth vector bundles. For each $f\in C^{\infty}(M)$ define
\[
\operatorname{ad}(f)\colon \operatorname{\mathbf{Op}}(\mathbf{E},\mathbf{F})\longrightarrow \operatorname{\mathbf{Op}}(\mathbf{E},\mathbf{F}),
\qquad
\operatorname{ad}(f)(P)=[P,f].
\]
Then $\operatorname{ad}(f)$ satisfies the following properties:
\begin{enumerate}[label=(\alph*)]
 \item $\mathbb{K}$-linearity:
 \[
 \operatorname{ad}(f)(\alpha P+\beta Q)
 =\alpha\,\operatorname{ad}(f)(P)
 +\beta\,\operatorname{ad}(f)(Q).
 \]

 \item $C^{\infty}(M)$-linearity in the operator:
 \[
 \operatorname{ad}(f)(hP+Q)
 =h\,\operatorname{ad}(f)(P)
 +\operatorname{ad}(f)(Q).
 \]

 \item $\mathbb{K}$-linearity in the function:
 \[
 \operatorname{ad}(\alpha f+\beta g)
 =\alpha\,\operatorname{ad}(f)
 +\beta\,\operatorname{ad}(g).
 \]

 \item If $c$ is constant, then
 \[
 \operatorname{ad}(c)=0.
 \]
\end{enumerate}
\end{proposition}

\begin{proof}
(a) This follows from property (a) of Proposition~\ref{prop:propiedades-conmutador}:
\[
\operatorname{ad}(f)(\alpha P+\beta Q)
=[\alpha P+\beta Q,f]
=\alpha [P,f]+\beta [Q,f].
\]

(b) Using property (b) of Proposition~\ref{prop:propiedades-conmutador},
\[
\operatorname{ad}(f)(hP+Q)
=[hP+Q,f]
=[hP,f]+[Q,f]
=h[P,f]+\operatorname{ad}(f)(Q).
\]

(c) Using property (c) of Proposition~\ref{prop:propiedades-conmutador},
\[
\operatorname{ad}(\alpha f+\beta g)(P)
=[P,\alpha f+\beta g]
=\alpha[P,f]+\beta[P,g].
\]

(d) If $c$ is constant, property (e) of Proposition~\ref{prop:propiedades-conmutador} implies
\[
\operatorname{ad}(c)(P)=[P,c]=0,
\]
for every $P\in \operatorname{\mathbf{Op}}(\mathbf{E},\mathbf{F})$.
\end{proof}

\begin{remark}\label{obs:operadores-diferenciales-en-haces-pdo}
 The preceding propositions imply that $\mathbf{PDO}^{(m)}(\mathbf{E},\mathbf{F})$ has the structure of an $C^{\infty}(M)$-module.
\end{remark}

We give some simple examples of partial differential operators in the Euclidean case to familiarize the reader with this concept.
\begin{example}\label{ejemplos de pdos sobre haz trivial}
Consider \(M=\mathbb{R}^n\) and the trivial bundles \(\mathbf{E}=\mathbf{F}=\underline{\mathbb{R}}\) of rank \(1\).

\begin{enumerate}
 \item The operator of multiplication by a smooth function is the most elementary example of a partial differential operator. Let \(M_{f}\) be defined by
 \[
 (M_{f}u)(x)=f(x)u(x),\qquad f\in C^\infty(\mathbb{R}^n).
 \]
 For any \(g,u\in C^\infty(\mathbb{R}^n)\) we have
 \[
 [M_{f},g](u)=M_{f}(gu)-g(M_{f}u)=f\cdot(gu)-g\cdot(fu)=0.
 \]
 Consequently, \([M_{f},g]=0\) for every \(g\in C^\infty(\mathbb{R}^n)\), and therefore
 \[
 M_{f}\in\mathbf{PDO}^{(0)}(\mathbf{E},\mathbf{F}).
 \]

 \item Now consider the partial derivative operator \(P=\partial_i=\displaystyle\frac{\partial}{\partial x_i}\). For \(u,f\in C^\infty(\mathbb{R}^n)\),
 \[
 [\partial_i,f](u)=\partial_i(fu)-f \partial_i u.
 \]
 Applying the product rule,
 \[
 \partial_i(fu)=(\partial_i f)u+f \partial_i u,
 \]
 gives
 \[
 [\partial_i,f](u)=(\partial_i f)u.
 \]
 That is, the commutator with \(f\) is simply multiplication by the function \(\partial_i f\), an operator of order \(0\). Thus, \([ \partial_i,f ]\in \mathbf{PDO}^{(0)}(\mathbf{E},\mathbf{F})\), and by definition
 \[
 \partial_i\in\mathbf{PDO}^{(1)}(\mathbf{E},\mathbf{F}).
 \]

 \item Finally, we analyze second-order partial derivatives. Let \(P=\partial_i^2\). For \(u\in C^\infty(\mathbb{R}^n)\),
 \[
 [\partial_i^2,f](u)=\partial_i^2(fu)-f \partial_i^2 u.
 \]
 Applying the product rule twice:
 \[
 \partial_i^2(fu)
 =\partial_i\big((\partial_i f)u+f \partial_i u\big)
 =(\partial_i^2 f)u+2(\partial_i f)\partial_i u+f \partial_i^2 u.
 \]
 Subtracting \(f \partial_i^2u\) gives
 \[
 [\partial_i^2,f](u)=(\partial_i^2 f)u+2(\partial_i f)\partial_i u.
 \]

 Taking a second commutator with another function \(g\in C^\infty(\mathbb{R}^n)\),
 \[
 [[\partial_i^2,f],g](u)
 =[M_{\partial_i^2 f}+2(\partial_i f)\partial_i, g](u).
 \]
 The first term \((\partial_i^2 f) \mathrm{Id}\) commutes with \(g\), since it is multiplication by a smooth function. For the second,
 \[
 [ 2(\partial_i f)\partial_i, g ](u)
 =2(\partial_i f)\big(\partial_i(gu)-g \partial_i u\big)
 =2(\partial_i f)(\partial_i g)u,
 \]
 which is multiplication by a smooth function, that is, an operator of order \(0\). Thus,

 $[[\partial_i^2,f],g]\in \mathbf{PDO}^{(0)}(\mathbf{E},\mathbf{F})$,
 y $\partial_i^2\in\mathbf{PDO}^{(2)}(\mathbf{E},\mathbf{F}).$

\end{enumerate}
\end{example}
The example of multiplication by a smooth function generalizes as follows:
\begin{example}\label{ejemplo multiplicacion es PDO0}
Let \(M\) be a smooth manifold with or without boundary and let \(\mathbf{E}\to M\) be a smooth vector bundle (real or complex). Fixing a function \(f\in C^\infty(M)\), consider the operator
\[
M_{f}\colon \Gamma(\mathbf{E})\longrightarrow \Gamma(\mathbf{E}),\qquad M_{f}(\mathbf{u}):=f \mathbf{u}.
\]
We wish to show that \(M_{f}\in \mathbf{PDO}^{(0)}(\mathbf{E},\mathbf{E})\). For every function \(g\in C^\infty(M)\) and every section \(\mathbf{u}\in \Gamma(\mathbf{E})\) we have
\[
[M_{f},g](\mathbf{u})=M_{f}(g\mathbf{u})-gM_{f}(\mathbf{u})=f(g\mathbf{u})-g(f\mathbf{u})=(fg)\mathbf{u}-(gf)\mathbf{u}=(fg)\mathbf{u}-(fg)\mathbf{u}=0.
\]
Therefore, \([M_{f},g]=0\) for every \(g\in C^\infty(M)\), meaning that
\[
M_{f}\in
\mathbf{PDO}^{(0)}(\mathbf{E},\mathbf{E}).
\]
\end{example}

Another example of a partial differential operator is the Laplace--Beltrami operator on a Riemannian manifold:

\begin{example}\label{ejemplo laplace beltrami pdo}
Let \((M,\mathbf{g})\) be a Riemannian manifold with or without boundary and consider the trivial bundle \(\mathbf{E}=\mathbf{F}=\underline{\mathbb{R}}\). The Laplace–Beltrami operator is defined by
\[
\Delta_{\mathbf{g}} u := \operatorname{div}(\operatorname{grad} u), \qquad u\in C^\infty(M).
\]
We will show that \(\Delta_{\mathbf{g}}\in \mathbf{PDO}^{(2)}(\mathbf{E},\mathbf{F})\) but that \(\Delta_{\mathbf{g}}\notin \mathbf{PDO}^{(1)}(\mathbf{E},\mathbf{F})\).

For every function \(f\in C^\infty(M)\), Proposition~\ref{leibniz para el laplaciano} gives
\[
\Delta_{\mathbf{g}}(fu)
=f \Delta_{\mathbf{g}} u + 2\langle\operatorname{grad}f,\operatorname{grad}u\rangle_{\mathbf{g}} + u \Delta_{\mathbf{g}} f,
\]
and therefore
\[
[\Delta_{\mathbf{g}},f](u)
=\Delta_{\mathbf{g}}(fu)-f \Delta_{\mathbf{g}} u
=2\langle\operatorname{grad}f,\operatorname{grad}u\rangle_{\mathbf{g}}+(\Delta_{\mathbf{g}} f) u.
\]
Now let $h\in C^{\infty}(M)$. By definition of the commutator, for every $u\in C^{\infty}(M)$ we have
\[
[[\Delta_{\mathbf{g}},f],h](u)
=
[\Delta_{\mathbf{g}},f](hu)
-
h[\Delta_{\mathbf{g}},f](u).
\]
Using the identity
\[
[\Delta_{\mathbf{g}},f](v)
=
2\langle\operatorname{grad}f,\operatorname{grad}v\rangle_{\mathbf{g}}
+
(\Delta_{\mathbf{g}}f)\,v,
\qquad v\in C^{\infty}(M),
\]
we obtain
\[
[\Delta_{\mathbf{g}},f](h u)
=
2\langle\operatorname{grad}f,\operatorname{grad}(hu)\rangle_{\mathbf{g}}
+
(\Delta_{\mathbf{g}}f)hu
\]
and
\[
[\Delta_{\mathbf{g}},f](u)
=
2\langle\operatorname{grad}f,\operatorname{grad}u\rangle_{\mathbf{g}}
+
(\Delta_{\mathbf{g}}f)\,u.
\]
Therefore,
\[
[[\Delta_{\mathbf{g}},f],h](u)
=
\Big(
2\langle\operatorname{grad}f,\operatorname{grad}(h u)\rangle_{\mathbf{g}}
+
(\Delta_{\mathbf{g}}f)h u
\Big)
-
h\Big(
2\langle\operatorname{grad}f,\operatorname{grad}u\rangle_{\mathbf{g}}
+
(\Delta_{\mathbf{g}}f)u
\Big)
\]
\[
=
2\langle\operatorname{grad}f,\operatorname{grad}(h u)\rangle_{\mathbf{g}}
+
(\Delta_{\mathbf{g}}f)h u
-
2h\langle\operatorname{grad}f,\operatorname{grad}u\rangle_{\mathbf{g}}
-
h(\Delta_{\mathbf{g}}f)\,u.
\]

We now use the Leibniz rule for the gradient (Proposition~\ref{leibniz gradiente}):
\[
\operatorname{grad}(h u)
=
h\operatorname{grad}u
+
u\,\operatorname{grad}h,
\]
so that
\[
\langle\operatorname{grad}f,\operatorname{grad}(hu)\rangle_{\mathbf{g}}
=
\big\langle\operatorname{grad}f,\,
 h\operatorname{grad}u
 +
 u\,\operatorname{grad}h
 \big\rangle_{\mathbf{g}}
\]
\[
=
h\langle\operatorname{grad}f,\operatorname{grad}u\rangle_{\mathbf{g}}
+
u\langle\operatorname{grad}f,\operatorname{grad}h\rangle_{\mathbf{g}},
\]
where we have used the $C^{\infty}(M)$--bilinearity of the inner product induced by $\mathbf{g}$. Consequently,
\[
2\langle\operatorname{grad}f,\operatorname{grad}(h u)\rangle_{\mathbf{g}}
=
2h\,\langle\operatorname{grad}f,\operatorname{grad}u\rangle_{\mathbf{g}}
+
2u\,\langle\operatorname{grad}f,\operatorname{grad}h\rangle_{\mathbf{g}}.
\]

Substituting this expression into the formula for $[[\Delta_{\mathbf{g}},f],h](u)$, we have
\[
[[\Delta_{\mathbf{g}},f],h](u)
=
\Big(
2h\langle\operatorname{grad}f,\operatorname{grad}u\rangle_{\mathbf{g}}
+
2u\,\langle\operatorname{grad}f,\operatorname{grad}h\rangle_{\mathbf{g}}
+
(\Delta_{\mathbf{g}}f)hu
\Big)
\]
\[
\qquad
-
\Big(
2h\langle\operatorname{grad}f,\operatorname{grad}u\rangle_{\mathbf{g}}
+
h(\Delta_{\mathbf{g}}f)\,u
\Big).
\]

We now simplify term by term. First,
\[
2h\langle\operatorname{grad}f,\operatorname{grad}u\rangle_{\mathbf{g}}
-
2h\langle\operatorname{grad}f,\operatorname{grad}u\rangle_{\mathbf{g}}
=0.
\]
Then, since multiplication of functions is commutative,
\[
(\Delta_{\mathbf{g}}f)h u
-
h(\Delta_{\mathbf{g}}f)\,u
=0.
\]
Therefore, only
\[
[[\Delta_{\mathbf{g}},f],h](u)
=
2u\,\langle\operatorname{grad}f,\operatorname{grad}h\rangle_{\mathbf{g}}.
\]
remains.

In other words, for every $u\in C^{\infty}(M)$ we have
\[
[[\Delta_{\mathbf{g}},f],h](u)
=
2\langle\operatorname{grad}f,\operatorname{grad}h\rangle_{\mathbf{g}}\,u,
\]
that is,
\[
[[\Delta_{\mathbf{g}},f],h]
=
M_{2\langle\operatorname{grad}f,\operatorname{grad}h\rangle_{\mathbf{g}}}.
\]

This result is multiplication by the smooth function \(2\langle\operatorname{grad}f,\operatorname{grad}h\rangle_{\mathbf{g}}\), so
\[
[[\Delta_{\mathbf{g}},f],h]\in \mathbf{PDO}^{(0)}(\mathbf{E},\mathbf{F})\quad \text{for any }f,h\in C^\infty(M).
\]
Thus,
\[
[\Delta_{\mathbf{g}},f]\in \mathbf{PDO}^{(1)}(\mathbf{E},\mathbf{F})\quad\text{for every }f\in C^{\infty}(M),
\]
and consequently
\[
\Delta_{\mathbf{g}}\in \mathbf{PDO}^{(2)}(\mathbf{E},\mathbf{F}).
\]

To show that \(\Delta_{\mathbf{g}}\notin \mathbf{PDO}^{(1)}(\mathbf{E},\mathbf{F})\), it suffices to note that if an operator belongs to \(\mathbf{PDO}^{(1)}(\mathbf{E},\mathbf{F})\), then its second commutators with smooth functions would vanish. However, we have obtained
\[
[[\Delta_{\mathbf{g}},f],h]=M_{\,2\langle\operatorname{grad}f,\operatorname{grad}h\rangle_{\mathbf{g}}},
\]
which is not the zero operator in general. Consequently,
\[
\Delta_{\mathbf{g}}\notin \mathbf{PDO}^{(1)}(\mathbf{E},\mathbf{F}).
\]
\end{example}
The preceding example shows that calling these partial differential operators of order less than or equal to $m$ makes sense, since having one of order less than or equal to $m$ does not necessarily imply that it has order less than or equal to $m-1$. Moreover, if it has order less than or equal to $m$, one can prove that it has order less than or equal to $m+1$, as will be done later in Lemma~\ref{lema:pdo monotonia}.

The covariant derivative defined on a smooth vector bundle is another example of a partial differential operator, in this case of order $\leq 1$:
\begin{proposition}\label{prop: conexion es pdo de orden 1}
 Let $M$ be a smooth manifold with or without boundary and let $\mathbf{E}\longrightarrow M$ be a smooth vector bundle over $M$. If \[\nabla\colon \Gamma(\mathbf{E})\longrightarrow \Gamma(T^{(0,1)}(TM)\otimes \mathbf{E})\] is a connection on $\mathbf{E}$, then $\nabla\in \mathbf{PDO}^{(1)}(\mathbf{E},T^{(0,1)}(TM)\otimes \mathbf{E})$.
\end{proposition}
\begin{proof}
 Let $f\in C^{\infty}(M)$ and $\mathbf{u}\in \Gamma(\mathbf{E})$. Then by Proposition~\ref{conexion equivalente}, \[[\nabla,f](\mathbf{u})=\nabla(f\mathbf{u})-f\nabla \mathbf{u}=df\otimes \mathbf{u}+f\nabla \mathbf{u}-f\nabla \mathbf{u}=df\otimes \mathbf{u}.\]

 Define $M_{df}^{\otimes}\in \mathbf{Op}(\mathbf{E},T^{(0,1)}(TM)\otimes \mathbf{E})$ by $M_{df}^{\otimes}(\mathbf{u})=df\otimes \mathbf{u}.$ Let $g\in C^{\infty}(M)$; then \[[M_{df}^{\otimes},g](\mathbf{u})=M_{df}^{\otimes}(g\mathbf{u})-gM_{df}^{\otimes}(\mathbf{u})=df\otimes (g\mathbf{u})-g(df\otimes \mathbf{u})=g(df\otimes \mathbf{u})-g(df\otimes \mathbf{u})=0,\], so $M_{df}^{\otimes}\in \mathbf{PDO}^{(0)}(\mathbf{E},T^{(0,1)}(TM)\otimes \mathbf{E})$ and consequently $\nabla\in\mathbf{PDO}^{(1)}(\mathbf{E},T^{(0,1)}(TM)\otimes \mathbf{E})$.
\end{proof}
To become more familiar with this concept, we study what partial differential operators of order $1$ look like on $\underline{\mathbb{R}}$ viewed as a vector bundle over $\mathbb{R}^{n}$.

\begin{example}\label{ej:operadores-diferenciales-en-haces-pdo}
 Consider again the trivial bundle $\underline{\mathbb{R}}$ over $\mathbb{R}^{n}$. As we saw earlier, the sections of $\underline{\mathbb{R}}$ are precisely the functions in $C^{\infty}(\mathbb{R}^{n})$. We wish to find $\mathbf{PDO}^{(1)}(\underline{\mathbb{R}},\underline{\mathbb{R}})$. Let $P\in \mathbf{PDO}^{(1)}(\underline{\mathbb{R}},\underline{\mathbb{R}})$ and $u,f\in C^{\infty}(\mathbb{R}^{n})$. Then \[[P,f]u=\sigma(f)u\quad\text{for some $\sigma(f)\in C^{\infty}(\mathbb{R}^{n})$},\], since $[P,f]$ is $C^{\infty}(M)$-linear because $P\in \mathbf{PDO}^{(1)}(\underline{\mathbb{R}},\underline{\mathbb{R}})$, and we may take $\sigma(f)=[P,f](1)$. On the other hand, for any $f,g\in C^{\infty}(\mathbb{R}^{n})$, \[\sigma(fg)u=[P,fg]u=[P,f](gu)+f([P,g]u)=\sigma(f)gu+f\sigma(g)u.\]

 Therefore, $\sigma(fg)=\sigma(f)g+f\sigma(g)$, so the map $\sigma$ is a derivation on $C^{\infty}(M)$; hence there exists a smooth vector field $\mathbf{X}\in \mathfrak{X}(\mathbb{R}^{n})$ such that \[\sigma(f)=\mathbf{X}\cdot f,\quad \forall f\in C^{\infty}(\mathbb{R}^{n}).\]

 Then for every $u\in C^{\infty}(\mathbb{R}^{n})$, \[P(u)=P(u\cdot 1)=[P,u]\cdot 1+u\cdot P(1)=\mathbf{X}\cdot u+P(1)\cdot u,\] where $P(1)\in C^{\infty}(\mathbb{R}^{n})$.
\end{example}

Operators of order $0$ can be characterized as follows:
\begin{lemma}\label{pdo y seccion}
Let $M$ be a smooth manifold with or without boundary and let $\mathbf{E},\mathbf{F}\longrightarrow M$ be smooth vector bundles. If $P\in \mathbf{PDO}^{(0)}(\mathbf{E},\mathbf{F})$, then there exists a section $\mathbf{A}\in \Gamma(\mathrm{Hom}(\mathbf{E},\mathbf{F}))$ such that
\[
P(\mathbf{u})(p)=\mathbf{A}(p)(\mathbf{u}(p)), \qquad \forall \mathbf{u}\in \Gamma(\mathbf{E}), p\in M.
\]
\end{lemma}
\begin{proof}
 By definition of a partial differential operator of order $0$, the map $P\colon \Gamma(\mathbf E)\longrightarrow\Gamma(\mathbf F)$ is $C^\infty(M)$-linear. Proposition~\ref{iso Hom-Gamma-Hom} identifies maps with this property with sections of $\operatorname{Hom}(\mathbf E,\mathbf F)$: there is a unique section $\mathbf A$ such that $P(\mathbf u)=\mathbf A(\mathbf u)$. Evaluating this equality at $p\in M$ gives $P(\mathbf u)(p)=\mathbf A(p)(\mathbf u(p))$.
\end{proof}
As a consequence of this lemma, partial differential operators are local.
\begin{proposition}\label{pdo es local}
 Let $M$ be a smooth manifold with or without boundary and let $\mathbf{E},\mathbf{F}\longrightarrow M$ be smooth vector bundles. Every $L\in\mathbf{PDO}^{(m)}(\mathbf{E},\mathbf{F})$ is a local operator, that is, $\forall \mathbf{u}\in \Gamma(\mathbf{E})$, \[\supp(L\mathbf{u})\subseteq\supp(\mathbf{u}).\]
\end{proposition}
 \begin{proof}
 The case $m=0$ follows from Lemma~\ref{pdo y seccion}, since if $\mathbf{u}(p)=0$ then $L(\mathbf{u})(p)=\mathbf{A}(p)(\mathbf{u}(p))=\mathbf{A}(p)(0)=0$, so \[\supp(L\mathbf{u})\subseteq\supp(\mathbf{u}).\]
 Suppose the assertion holds for $m\geq 0$ and let $L\in \mathbf{PDO}^{(m+1)}(\mathbf{E},\mathbf{F})$.

 For any $f\in C^{\infty}(M)$ we have \[L(f\mathbf{u})=[L,f]\mathbf{u}+fL(\mathbf{u}).\]

 Since $[L,f]\in \mathbf{PDO}^{(m)}(\mathbf{E},\mathbf{F})$, the induction hypothesis gives
 \begin{equation}\label{eq: contencion soportes localidad pdo}
 \begin{aligned}
 \supp L(f\mathbf{u})&\subseteq \supp([L,f]\mathbf{u})\cup \supp(fL(\mathbf{u}))\subseteq \supp(\mathbf{u})\cup(\supp(f)\cap\supp(L\mathbf{u}))\\
 &\subseteq \supp(\mathbf{u})\cup \supp(f),\quad \forall f\in C^{\infty}(M).
 \end{aligned}
 \end{equation}

Under these hypotheses, we show that $\operatorname{supp}(L\mathbf{u})\subseteq \operatorname{supp}(\mathbf{u})$. We argue by contradiction and suppose that there exists $p\in M$ such that $p\in \supp(L\mathbf{u})$ but $p\notin \supp(\mathbf{u})$. Since $\supp(\mathbf{u})$ is closed and $M$ is a regular topological space (in fact, it is metrizable), there is an open set $U\subseteq M$ such that $p\in U$ and $\overline{U}\cap \supp(\mathbf{u})=\varnothing$, that is, $\supp(\mathbf{u})\subseteq M\setminus \overline{U}$, where $M\setminus \overline{U}$ is open. By Proposition~\ref{funciones flan}, there exists $f\in C^{\infty}(M)$ such that $f\equiv 0$ on $\overline{U}$ and $f\equiv 1$ on $\supp(\mathbf{u})$. The identity $f\mathbf{u}=\mathbf{u}$ holds throughout $M$, since $f\equiv 1$ on $\supp(\mathbf{u})$, and if a point $x$ lies outside $\supp(\mathbf{u})$, then $\mathbf{u}(x)=0=f(x)\mathbf{u}(x)$. Therefore, \[\supp(L\mathbf{u})=\supp(L(f\mathbf{u})),\], allowing us to conclude that $p\in \supp(L(f\mathbf{u}))$ since we had taken $p\in \supp(L\mathbf{u})$. But this is impossible, because \eqref{eq: contencion soportes localidad pdo} would imply $p\in \supp(\mathbf{u})\cup \supp(f)$, whereas $p\notin \supp(\mathbf{u})$ and $p\notin \supp(f)$ (the latter holds because $f\equiv 0$ on $U\subseteq \overline{U}$, so $U\cap \{x\in M\mid f(x)\neq 0\}=\varnothing$, and thus $p\notin\overline{\{x\in M\mid f(x)\neq 0\}}=\supp(f)$). The contradiction arose from assuming that $p\notin \supp(\mathbf{u})$, so $\supp(L\mathbf{u})\subseteq \supp(\mathbf{u})$, completing the proof.
 \end{proof}
\begin{remark}\label{remark: pdo es local implica igualdad de pdos si hay igualdad de funciones}
 If $L\in \mathbf{PDO}^{(m)}(\mathbf{E},\mathbf{F})$ and $\mathbf{u},\mathbf{v}\in \Gamma(\mathbf{E})$ are such that there is an open set $U\subseteq M$ with $\mathbf{u}\equiv \mathbf{v}$ on $U$, and we define $\mathbf{w}:=\mathbf{u}-\mathbf{v}$, this shows that $U\cap \supp(\mathbf{w})=\varnothing$. Proposition~\ref{pdo es local} therefore implies that $\supp(L\mathbf{w})\subseteq \supp(\mathbf{w})$, and consequently $U\cap \supp(L\mathbf{w})=\varnothing$, from which we conclude that $L\mathbf{u}-L\mathbf{v}=L\mathbf{w}=0$ on $U$; that is, $L\mathbf{u}$ and $L\mathbf{v}$ agree on $U$.

 This justifies the word ``local'' in that proposition, since it tells us that the values of a partial differential operator are completely determined by the values of the section to which it is applied on an open subset of $M$.
\end{remark}
We show that the composition of partial differential operators is again a partial differential operator. For this, we need the following lemmas:
\begin{lemma}\label{lema:conmutador-composicion}
 Let $M$ be a smooth manifold with or without boundary. Let $\mathbf{E},\mathbf{F},\mathbf{G}$ be smooth vector bundles over $M$, let $f\in C^\infty(M)$, and let $P\colon \Gamma(\mathbf{F})\longrightarrow\Gamma(\mathbf{G})$, $Q\colon \Gamma(\mathbf{E})\longrightarrow\Gamma(\mathbf{F})$ be linear operators. Then, for every $\mathbf{u}\in \Gamma(\mathbf{E})$,
\[
[P\circ Q,f](\mathbf{u})=[P,f]\big(Q(\mathbf{u})\big)+P\big([Q,f](\mathbf{u})\big).
\]
\end{lemma}

\begin{proof}
By definition,
\[
[P\circ Q,f](\mathbf{u})=(P\circ Q)(f\mathbf{u})-f(P\circ Q)(\mathbf{u})=P\big(Q(f\mathbf{u})\big)-fP\big(Q(\mathbf{u})\big).
\]
Add and subtract $P\big(fQ(\mathbf{u})\big)$:
\[
\begin{aligned}
[P\circ Q,f](\mathbf{u})
&= \Big(P\big(Q(f\mathbf{u})\big)-P\big(fQ(\mathbf{u})\big)\Big)+\Big(P\big(fQ(\mathbf{u})\big)-fP\big(Q(\mathbf{u})\big)\Big)\\
&= P\big(Q(f\mathbf{u})-fQ(\mathbf{u})\big)+\big(P(fQ(\mathbf{u}))-fP(Q(\mathbf{u}))\big)\\
&= P\big([Q,f](\mathbf{u})\big)+[P,f]\big(Q(\mathbf{u})\big).
\end{aligned}
\]
\end{proof}
Also note that if $P\in \mathbf{PDO}^{(m)}(\mathbf{E},\mathbf{F})$, this should be interpreted as saying that $P$ is a partial differential operator of order less than or equal to $m$. This tells us that if a partial differential operator has order less than or equal to $m$, then it has order less than or equal to $m+1$, justifying the initial choice of terminology.
\begin{lemma}\label{lema:pdo monotonia}
Let $M$ be a smooth manifold with or without boundary, let $\mathbf{E},\mathbf{F}\longrightarrow M$ be smooth vector bundles, and let $m\geq 0$. Then \[\mathbf{PDO}^{(m)}(\mathbf{E},\mathbf{F})\subseteq \mathbf{PDO}^{(m+1)}(\mathbf{E},\mathbf{F}).\]
\end{lemma}

\begin{proof}
We prove the assertion by induction on $m$. If $m=0$ and $P\in \mathbf{PDO}^{(0)}(\mathbf{E},\mathbf{F})$, then for every $f\in C^{\infty}(M)$ $[P,f]=0$, so by the recursive definition of partial differential operators we have $[P,f]=0\in \mathbf{PDO}^{(0)}(\mathbf{E},\mathbf{F})$. Suppose as the induction hypothesis that $\mathbf{PDO}^{(m)}(\mathbf{E},\mathbf{F})\subseteq \mathbf{PDO}^{(m+1)}(\mathbf{E},\mathbf{F})$. If $P\in \mathbf{PDO}^{(m+1)}(\mathbf{E},\mathbf{F})$, then for every $f\in C^{\infty}(M)$, $[P,f]\in \mathbf{PDO}^{(m)}(\mathbf{E},\mathbf{F})\subseteq \mathbf{PDO}^{(m+1)}(\mathbf{E},\mathbf{F})$, so $[P,f]\in \mathbf{PDO}^{(m+1)}(\mathbf{E},\mathbf{F})$ $\forall f\in C^{\infty}(M)$. This allows us to conclude that $P\in \mathbf{PDO}^{(m+2)}(\mathbf{E},\mathbf{F})$.
\end{proof}
We now examine their behavior under composition.
\begin{proposition}\label{prop: composicion de pdos es pdo}
 Let $M$ be a smooth manifold with or without boundary and let $\mathbf{E},\mathbf{F},\mathbf{G}$ be smooth real or complex vector bundles over $M$. If $P\in \mathbf{PDO}^{(m)}(\mathbf{F},\mathbf{G})$ and $Q\in \mathbf{PDO}^{(n)}(\mathbf{E},\mathbf{F})$, then $P\circ Q\in \mathbf{PDO}^{(m+n)}(\mathbf{E},\mathbf{G})$.
\end{proposition}
\begin{proof}
Let $k=m+n$. We proceed by induction on $k$.

For the base case $k=0$, we have $m=n=0$, so $P,Q$ are $C^{\infty}(M)$–linear, and their composition is as well. Therefore, $P\circ Q\in \mathbf{PDO}^{(0)}(\mathbf{E},\mathbf{G})$.

Now suppose as the induction hypothesis that the result holds for any $m,n\geq 0$ such that $m+n=k$; that is, if $P\in \mathbf{PDO}^{(m)}(\mathbf{F},\mathbf{G})$ and $Q\in \mathbf{PDO}^{(n)}(\mathbf{E},\mathbf{F})$ with $m+n=k$, then $P\circ Q\in \mathbf{PDO}^{(k)}(\mathbf{E},\mathbf{G})$. We must prove that the result holds for any $m,n\geq 0$ with $m+n=k+1$.

Let $P\in \mathbf{PDO}^{(m)}(\mathbf{F},\mathbf{G})$ and $Q\in \mathbf{PDO}^{(n)}(\mathbf{E},\mathbf{F})$ with $m+n=k+1$. For $f\in C^{\infty}(M)$, Lemma~\ref{lema:conmutador-composicion} states that
\[
[P\circ Q,f]=[P,f]\circ Q+P\circ[Q,f].
\]
If $m=0$, the $C^\infty(M)$--linearity of $P$ gives $[P,f]=0$; if $m\geq1$, then $[P,f]\in\mathbf{PDO}^{(m-1)}(\mathbf{F},\mathbf{G})$ and $(m-1)+n=k$, so the induction hypothesis gives $[P,f]\circ Q\in\mathbf{PDO}^{(k)}(\mathbf{E},\mathbf{G})$. For the second summand, if $n=0$ then $[Q,f]=0$; if $n\geq1$, then $[Q,f]\in\mathbf{PDO}^{(n-1)}(\mathbf{E},\mathbf{F})$ and $m+(n-1)=k$, so $P\circ[Q,f]\in\mathbf{PDO}^{(k)}(\mathbf{E},\mathbf{G})$ by the same hypothesis.

Consequently, $[P\circ Q,f]\in \mathbf{PDO}^{(k)}(\mathbf{E},\mathbf{G})$ for every $f\in C^{\infty}(M)$. The recursive definition implies that $P\circ Q\in \mathbf{PDO}^{(k+1)}(\mathbf{E},\mathbf{G})=\mathbf{PDO}^{(m+n)}(\mathbf{E},\mathbf{G})$.
\end{proof}
As a corollary of the preceding results, the most common differential operators on $\mathbb{R}^{n}$ defined in terms of partial derivatives are partial differential operators in the sense we have defined.
\begin{corollary}\label{derivada parcial alpha esima es pdo}
 The operator \[L=\displaystyle\sum_{|\alpha|\leq m}a_{\alpha}(x)D^{\alpha}\colon C^{\infty}(\mathbb{R}^{n})\longrightarrow C^{\infty}(\mathbb{R}^{n})\] with $a_\alpha \in C^\infty (\mathbb{R}^n), $ is a partial differential operator of order less than or equal to $m$.
\end{corollary}
\begin{proof}
 By Example~\ref{ejemplos de pdos sobre haz trivial}, each partial derivative $\partial_{i}$ is a first-order partial differential operator, so Proposition~\ref{prop: composicion de pdos es pdo} ensures that $L$ is a partial differential operator of order $\leq m$, since multiplication by a smooth function does not change the order of $\mathbf{PDO}$.
\end{proof}
We examine what happens when the ad operators are permuted:
\begin{lemma}\label{lem:operadores-diferenciales-en-haces-pdo}
 Let $M$ be a smooth manifold with or without boundary and let $\mathbf{E},\mathbf{F}\longrightarrow M$ be smooth vector bundles over $M$. Then for any $P\in \mathbf{PDO}(\mathbf{E},\mathbf{F})$ and any $f,g\in C^{\infty}(M)$ \[(\operatorname{ad}(f)\circ\operatorname{ad}(g))P=(\operatorname{ad}(g)\circ\operatorname{ad}(f))P.\]
\end{lemma}
\begin{proof} Let $P\in \mathbf{PDO}(\mathbf{E},\mathbf{F})$, $f,g\in C^{\infty}(M)$, and $\mathbf{u}\in \Gamma(\mathbf{E})$. Then
\[(\operatorname{ad}(f)\circ \operatorname{ad}(g))(P)(\mathbf{u})=[[P,g],f](\mathbf{u})=[P,g](f\mathbf{u})-f[P,g](\mathbf{u})=\color{red}{P(gfu)}-\color{blue}{gP(f\mathbf{u})}-\color{violet}{fP(g\mathbf{u})}+\color{teal}{fgP(\mathbf{u})}.\]

On the other hand, \[(\operatorname{ad}(g)\circ \operatorname{ad}(f))(P)(\mathbf{u})=[[P,f],g](\mathbf{u})=[P,f](g\mathbf{u})-g[P,f](\mathbf{u})=\color{red}{P(fgu)}-\color{violet}{fP(g\mathbf{u})}-\color{blue}{gP(f\mathbf{u})}+\color{teal}{gfP(\mathbf{u})}.\]
Terms of the same color are equal, since we can interchange $f$ and $g$, so $(\operatorname{ad}(f)\circ\operatorname{ad}(g))P=(\operatorname{ad}(g)\circ\operatorname{ad}(f))P$.

\end{proof}
\begin{remark}\label{obs: ad es invariante bajo permutaciones de f}
 The preceding lemma tells us that if $P\in \mathbf{PDO}^{(m)}(\mathbf{E},\mathbf{F})$, then for any $f_{1},\dots,f_{m}\in C^{\infty}(M)$, the bundle homomorphism \[\operatorname{ad}(f_{1})\circ\operatorname{ad}(f_{2})\cdots\circ\operatorname{ad}(f_{m})P\] is unchanged when $f_{1},\dots,f_{m}$ are permuted.
\end{remark}
We next characterize iterated commutators using equality of the respective differentials of the functions; to do so, we first prove a lemma.
\begin{lemma}\label{lema: ideales pdo}
 Let $M$ be a smooth manifold with or without boundary. For each $p\in M$, consider the ideals of the ring $C_{p}^{\infty}(M)$,
 \[
 \mathbf{m}_{p}
 =
 \{[f]\in C_{p}^{\infty}(M) \mid f(p)=0\},
 \qquad
 \mathcal{I}_{p}
 :=
 \{[f]\in C_{p}^{\infty}(M)\mid f(p)=0,\; df(p)=0\}.
 \]
 Then $\mathcal{I}_{p}=\mathbf{m}_{p}^{2}$; that is, every smooth function $f$ that vanishes at $p$ together with its differential at $p$ can be written, in a neighborhood of $p$, as
 \[
 f=\displaystyle\sum_{j=1}^{k} g_{j}h_{j},
 \qquad
 g_{j},h_{j}\in \mathbf{m}_{p}.
 \]
\end{lemma}

\begin{proof}
Let $p\in M$. First observe that $\mathbf{m}_{p}^{2}\subseteq \mathcal{I}_{p}$. If $g,h\in \mathbf{m}_{p}$, then $g(p)=h(p)=0$, and the product rule gives
\[
d(gh)(p)=g(p)\, dh(p)+h(p)\, dg(p)=0.
\]
Therefore, $gh\in \mathcal{I}_{p}$, and every finite combination of such products belongs to $\mathcal{I}_{p}$.

Now let $f\in \mathcal{I}_{p}$, so that $f(p)=0$ and $df(p)=0$. Take a smooth chart $(U,\phi)$ around $p$, chosen so that $\phi(U)$ is convex, $\phi(p)=0\in\mathbb{R}^{n}$, and write
\[
\phi=(x^{1},\dots,x^{n}).
\]
Consider the coordinate representation of $f$ given by
\[
\tilde{f}:=f\circ \phi^{-1}\colon \phi(U)\longrightarrow \mathbb{R}.
\]
The conditions $f(p)=0$ and $df(p)=0$ are equivalent to
\[
\tilde{f}(0)=0,\qquad
\frac{\partial \tilde{f}}{\partial x^{i}}(0)=0,\quad \forall i\in\{1,\dots,n\}.
\]

Taylor's theorem (Theorem~\ref{taylor multivariable multiindices}) for $k=1$, applied at the point $a=0$, gives
\[
\tilde{f}(x)
=\tilde{f}(0)+\displaystyle\sum_{i=1}^{n}\frac{\partial \tilde{f}}{\partial x^{i}}(0)\, x^{i}
+\displaystyle\sum_{|\alpha|=2}\frac{2}{\alpha !}\, x^{\alpha}
\int_{0}^{1}(1-t)\, D^{\alpha}\tilde{f}(t x)\, dt.
\]
The first two terms vanish, so, passing from multi-index notation to the usual Leibniz notation, we obtain
\[
\tilde{f}(x)
=\displaystyle\sum_{i,j=1}^{n}x^{i}x^{j}R_{ij}(x),
\]
where
\[
R_{ij}(x):=\int_{0}^{1}(1-t)\,
\frac{\partial^{2}\tilde{f}}{\partial x^{i}\partial x^{j}}(t x)\, dt.
\]
Each function $R_{ij}$ is of class $C^{\infty}$ in a neighborhood of $0$ because $\tilde{f}$ is.

Define
\[
\tilde{g}_{i}(x):=\displaystyle\sum_{j=1}^{n}x^{j}R_{ij}(x),
\]
so that
\[
\tilde{f}(x)=\displaystyle\sum_{i=1}^{n}x^{i}\tilde{g}_{i}(x).
\]
Observe that, for each $i$,
\[
\tilde{g}_{i}(0)=\displaystyle\sum_{j=1}^{n} 0\cdot R_{ij}(0)=0.
\]

Now let $g_{i}:=\tilde{g}_{i}\circ\phi\in C^{\infty}(U)$. Then, for every $q\in U$ we have
\[
f(q)=\displaystyle\sum_{i=1}^{n}x^{i}(q)\, g_{i}(q).
\] We can extend $x^{i},g_{i}$ smoothly to $M$ using Theorem~\ref{extension de funciones suaves} so that this equality holds in a neighborhood $W\subseteq U$ of $p$.

Therefore,
\[
[f]=\left[\displaystyle\sum_{i=1}^{n}x^{i}g_{i}\right]
=\displaystyle\sum_{i=1}^{n}[x^{i}]\,[g_{i}]
\]
as germs of functions. Since $x^{i}(p)=0$ and $g_{i}(p)=\tilde{g}_{i}(0)=0$, we have $[x^{i}],[g_{i}]\in \mathbf{m}_{p}$ for every $i\in \{1,\dots,n\}$, and therefore
\[
[f]\in \mathbf{m}_{p}^{2}.
\]

We have proved that $\mathcal{I}_{p}\subseteq \mathbf{m}_{p}^{2}$, and together with the reverse inclusion we obtain
\[
\mathcal{I}_{p}=\mathbf{m}_{p}^{2}.
\]
\end{proof}
\begin{definition}\label{def:operadores-diferenciales-en-haces-tal-cual-lema-lema-seccion-valor}\index{iterated commutator!pointwise evaluation}
 Let $M$ be a smooth manifold with or without boundary, let $\mathbf{E},\mathbf{F}\longrightarrow M$ be smooth vector bundles, and let $P\in \mathbf{PDO}^{(m)}(\mathbf{E},\mathbf{F})$. For $p\in M$ and $f_{1},\dots,f_{m}\in C^{\infty}(M)$ define
 \[
 \operatorname{ad}(f_{1})\cdots \operatorname{ad}(f_{m})(P)|_{p}\colon \mathbf{E}_{p}\longrightarrow \mathbf{F}_{p}
 \]
 by
 \[
 \bigl(\operatorname{ad}(f_{1})\cdots \operatorname{ad}(f_{m})(P)|_{p}\bigr)(v)
 :=\bigl(\operatorname{ad}(f_{1})\cdots \operatorname{ad}(f_{m})(P)\bigr)(\boldsymbol{\sigma})(p),
 \]
 where $\boldsymbol{\sigma}\in \Gamma(\mathbf{E})$ satisfies $\boldsymbol{\sigma}(p)=v$. Such a section exists by Lemma~\ref{lema:seccion-valor-prescrito}, and the value does not depend on the choice of $\boldsymbol{\sigma}$ because the iterated operator has order zero.
\end{definition}
\begin{proposition}\label{prop: localidad respecto al diferencial del adjunto de Lie}
 Let $M$ be a smooth manifold with or without boundary, let $\mathbf{E},\mathbf{F}\longrightarrow M$ be smooth vector bundles, and let $m\geq 1$. Let also $P\in\mathbf{PDO}^{(m)}(\mathbf{E},\mathbf{F})$, $f_{i},g_{i}\in C^{\infty}(M)$ with $i\in \{1,\dots,m\}$ be such that for $p\in M$, \[df_{i}(p)=dg_{i}(p)\in T_{p}^{*}M \qquad \forall i\in \{1,\dots,m\},\]
 then \[\operatorname{ad}(f_{1})\circ\operatorname{ad}(f_{2})\cdots \circ\operatorname{ad}(f_{m})P|_{p}=\operatorname{ad}(g_{1})\circ\operatorname{ad}(g_{2})\cdots \circ\operatorname{ad}(g_{m})P|_{p}.\]
\end{proposition}
\begin{proof}
By Proposition~\ref{prop:propiedades-ad}, $\operatorname{ad}(c)=0$ for every constant function $c$. We may replace each $g_i$ by
\[
 \widehat g_i:=g_i+f_i(p)-g_i(p).
\]
This replacement does not change $\operatorname{ad}(g_i)$ and simultaneously ensures
\[
 f_i(p)=\widehat g_i(p),
 \qquad df_i(p)=d\widehat g_i(p).
\]
We again write $g_i$ in place of $\widehat g_i$ and set $\varphi_i:=f_i-g_i$. Then
\[
 \varphi_i(p)=0,\qquad d\varphi_i(p)=0.
\]
By Lemma~\ref{lema: ideales pdo}, the germ of $\varphi_i$ belongs to $\mathbf m_p^2$.

Let $Q\in\mathbf{PDO}^{(1)}(\mathbf E,\mathbf F)$ and let $\varphi\in C^\infty(M)$ with $\varphi(p)=0$ and $d\varphi(p)=0$. The same lemma provides, in a neighborhood of $p$, functions $\alpha_\nu,\beta_\nu$ such that
\[
 \varphi=\sum_{\nu=1}^{q}\alpha_\nu\beta_\nu,
 \qquad \alpha_\nu(p)=\beta_\nu(p)=0.
\]
Theorem~\ref{extension de funciones suaves} allows us to choose global representatives of these germs. By Remark~\ref{remark: pdo es local implica igualdad de pdos si hay igualdad de funciones}, equality in a neighborhood of $p$ suffices to compute the value at $p$. The Leibniz rule in Proposition~\ref{prop:propiedades-conmutador} gives
\begin{equation}
 \operatorname{ad}(\varphi)(Q)
 =\sum_{\nu=1}^{q}
 \left(
 \operatorname{ad}(\alpha_\nu)(Q)\circ M_{\beta_\nu}
 +M_{\alpha_\nu}\circ\operatorname{ad}(\beta_\nu)(Q)
 \right),
 \label{eq:paso-sustitucion-conmutadores-diferenciales}
\end{equation}
where $M_h$ denotes multiplication by $h$. Since $Q$ has order at most one, the operators $\operatorname{ad}(\alpha_\nu)(Q)$ and $\operatorname{ad}(\beta_\nu)(Q)$ have order zero. If $\boldsymbol\sigma\in\Gamma(\mathbf E)$, evaluating \eqref{eq:paso-sustitucion-conmutadores-diferenciales} at $p$ gives
\[
 \bigl(\operatorname{ad}(\alpha_\nu)(Q)\bigr)
       (\beta_\nu\boldsymbol\sigma)(p)=0,
 \qquad
 \alpha_\nu(p)
 \bigl(\operatorname{ad}(\beta_\nu)(Q)\boldsymbol\sigma\bigr)(p)=0.
\]
Hence,
\[
 \operatorname{ad}(\varphi)(Q)|_p=0.
\]

For $r\in\{0,\dots,m\}$ define
\[
 T_r:=
 \operatorname{ad}(g_1)\cdots\operatorname{ad}(g_r)
 \operatorname{ad}(f_{r+1})\cdots\operatorname{ad}(f_m)(P),
\]
with the convention that an empty product is the identity, and for $r\in\{1,\dots,m\}$ set
\[
 Q_r:=
 \operatorname{ad}(g_1)\cdots\operatorname{ad}(g_{r-1})
 \operatorname{ad}(f_{r+1})\cdots\operatorname{ad}(f_m)(P).
\]
Here exactly $m-1$ commutators have been applied to an operator of order at most $m$, so $Q_r$ has order at most one. Moreover, the operators $\operatorname{ad}(h)$ commute with one another by Lemma~\ref{lem:operadores-diferenciales-en-haces-pdo}. Also using linearity from Proposition~\ref{prop:propiedades-ad},
\[
 T_{r-1}-T_r
 =\operatorname{ad}(f_r-g_r)(Q_r)
 =\operatorname{ad}(\varphi_r)(Q_r).
\]
Equality $\operatorname{ad}(\varphi)(Q)|_p=0$, applied to $\varphi=\varphi_r$ and $Q=Q_r$, shows that $(T_{r-1}-T_r)|_p=0$ for every $r$. Summing telescopically,
\[
 (T_0-T_m)|_p
 =\sum_{r=1}^{m}(T_{r-1}-T_r)|_p=0.
\]
Since $T_0$ and $T_m$ are, respectively, the two iterated commutators in the statement, the required equality is proved.

\end{proof}

We will characterize partial differential operators between the trivial vector bundles $\underline{\mathbb{K}}^{p}$ and $\underline{\mathbb{K}}^{q}$ explicitly in terms of their local coordinates, giving us a better idea of what a partial differential operator looks like in general.

The ideals in Lemma~\ref{lema: ideales pdo} and the use of Taylor's theorem in its proof suggest that these ideals are important in characterizing functions in terms of their derivatives. Just as the cotangent bundle captures information about the first derivatives of smooth functions, there are bundles called \textit{jet bundles} that store information about higher-order derivatives of a function. Introducing and studying these bundles will make it easier to characterize partial differential operators locally and to work with related concepts both intrinsically and locally.
\begin{definition}[Jets of order $m$ of smooth functions]\label{def:jet-funciones}\index{jet!of finite order}
Let $M$ be a smooth manifold with or without boundary and let $p\in M$. Denote the ring of germs of smooth functions at $p$ by $C_{p}^{\infty}(M)$, and write
\[
\mathbf{m}_{p}
=
\{[f]\in C_{p}^{\infty}(M)\mid f(p)=0\}
\]
which is a maximal ideal.

For $m\ge 0$ define the power
\[
\mathbf{m}_{p}^{\,m+1}
:=
\Big\langle
[f_{1}]\cdots [f_{m+1}]
\;\mid\;
[f_{i}]\in \mathbf{m}_{p},\ i=1,\dots,m+1
\Big\rangle
\subset C_{p}^{\infty}(M).
\]

The \textit{space of jets of order $m$ of smooth functions at $p$} is defined as the quotient vector space
\[
J^{m}_{p}(M)
=
C_{p}^{\infty}(M)\big/\mathbf{m}_{p}^{\,m+1}.
\]
For $[f]\in C_{p}^{\infty}(M)$, denote its class in $J^{m}_{p}(M)$ by
\[
j^{m}_{p}f,
\]
and call it the \textit{jet of order $m$} of $f$ at $p$.
\end{definition}
This definition can be characterized in terms of derivatives, showing the importance of jets.

\begin{lemma}\label{lema:igualdad-ideales-jets}
 Let $M$ be a smooth manifold with or without boundary. For each $p\in M$, consider the subsets of the ring $C_{p}^{\infty}(M)$
 \[
 \mathcal{I}_{p}^{(m)}
 :=
 \bigl\{[f]\in C_{p}^{\infty}(M)\;\mid\;
 D^{\alpha}f(p)=0\text{ for every }|\alpha|\le m\bigr\},
 \qquad
 \mathbf{m}_{p}^{\,m+1}
 :=
 \big\langle [f_{1}]\cdots[f_{m+1}]
 \;\mid\;
 [f_{i}]\in \mathbf{m}_{p}\big\rangle .
 \]
 Then $\mathcal{I}_{p}^{(m)}=\mathbf{m}_{p}^{\,m+1}$.
\end{lemma}

\begin{proof}
Let $p\in M$. First we show that $\mathbf{m}_{p}^{\,m+1}\subseteq \mathcal{I}_{p}^{(m)}$. If
\[
[f]=[g_{1}]\cdots[g_{m+1}], \qquad [g_{i}]\in \mathbf{m}_{p},
\]
then $g_{i}(p)=0$ for every $i\in \{1,\dots,m+1\}$. For any derivative of order $|\alpha|\le m$, the Leibniz rule in Proposition~\ref{prop:leibniz-multiindices-varios-factores} gives
\[
D^{\alpha}(g_{1}\cdots g_{m+1})(p)
=
\displaystyle\sum_{\alpha_{1}+\cdots+\alpha_{m+1}=\alpha}
\binom{\alpha}{\alpha_{1},\dots,\alpha_{m+1}}\,
D^{\alpha_{1}}g_{1}(p)\cdots D^{\alpha_{m+1}}g_{m+1}(p),
\]
where
\[
\binom{\alpha}{\alpha_{1},\dots,\alpha_{m+1}}
=
\frac{\alpha!}{\alpha_{1}!\cdots \alpha_{m+1}!}
\]
is the multinomial coefficient associated with the multi-indices $\alpha,\alpha_{1},\dots,\alpha_{m+1}$. Since $|\alpha|\le m$, in every tuple \(\alpha_{1},\dots,\alpha_{m+1}\) with sum \(\alpha\) at least one of the multi-indices is zero. This implies that in each term of the Leibniz rule, at least one factor $g_i(p)$ is left undifferentiated; since $g_i(p)=0$, all summands vanish. Therefore,
\[
D^{\alpha}(g_{1}\cdots g_{m+1})(p)=0,\qquad \forall\,|\alpha|\le m,
\]
which shows that $[f]\in \mathcal{I}_{p}^{(m)}$.

Now let $[f]\in \mathcal{I}_{p}^{(m)}$, so that
\[
D^{\alpha}f(p)=0,\qquad \forall |\alpha|\le m.
\]
Take a smooth chart $(U,\phi)$ around $p$, with $\phi(U)$ convex and $\phi(p)=0$, and write
\[
\phi=(x^{1},\dots,x^{n}).
\]
Let
\[
\tilde{f}:=f\circ \phi^{-1}\colon \phi(U)\longrightarrow\mathbb{R}.
\]
Then
\[
D^{\alpha}\tilde{f}(0)=0,\qquad \forall |\alpha|\le m.
\]

Taylor's theorem (Theorem~\ref{taylor multivariable multiindices}) for $k=m$ at the point $0$ gives
\[
\tilde{f}(x)
=
\displaystyle\sum_{|\alpha|\le m}
\frac{D^{\alpha}\tilde{f}(0)}{\alpha!}\,x^{\alpha}
+
\displaystyle\sum_{|\alpha|=m+1}
\frac{m+1}{\alpha!}\,x^{\alpha}
\int_{0}^{1}(1-t)^{m}\,D^{\alpha}\tilde{f}(t x)\,dt.
\]
The first summand vanishes by hypothesis, so
\[
\tilde{f}(x)
=
\displaystyle\sum_{|\alpha|=m+1}
x^{\alpha}\,\tilde{g}_{\alpha}(x),
\]
where
\[
\tilde{g}_{\alpha}(x)
:=
\frac{m+1}{\alpha!}
\int_{0}^{1}(1-t)^{m}\,D^{\alpha}\tilde{f}(t x)\,dt,
\]
and each $\tilde{g}_{\alpha}$ is of class $C^{\infty}$ in a neighborhood of $0$.

Let $g_{\alpha}:=\tilde{g}_{\alpha}\circ \phi\in C^{\infty}(U)$. Then, for every $q\in U$,
\[
f(q)=\displaystyle\sum_{|\alpha|=m+1}x^{\alpha}(q)\, g_{\alpha}(q).
\]
Using Theorem~\ref{extension de funciones suaves}, extend $x^{i}$ and $g_{\alpha}$ to global smooth functions without changing this equality in a neighborhood $W\subseteq U$ of $p$.

Thus,
\[
[f]
=
\left[\displaystyle\sum_{|\alpha|=m+1}x^{\alpha} g_{\alpha}\right]
=
\displaystyle\sum_{|\alpha|=m+1}[x^{\alpha}]\, [g_{\alpha}]
\quad\text{in }C_{p}^{\infty}(M).
\]

Observe that if $|\alpha|=m+1$, then
\[
x^{\alpha}=(x^{1})^{\alpha_{1}}\cdots(x^{n})^{\alpha_{n}}
\]
is a product of exactly $m+1$ coordinate functions, all in $\mathbf{m}_{p}$. Therefore,
\[
[x^{\alpha}] \in \mathbf{m}_{p}^{\,m+1}.
\]

Finally, since $\mathbf{m}_{p}^{\,m+1}$ is an ideal of $C_{p}^{\infty}(M)$, multiplying an element of $\mathbf{m}_{p}^{\,m+1}$ by any germ $[g_{\alpha}]$ again yields an element of $\mathbf{m}_{p}^{\,m+1}$. It follows that
\[
[f]\in \mathbf{m}_{p}^{\,m+1}.
\]

This proves that $\mathcal{I}_{p}^{(m)}\subseteq \mathbf{m}_{p}^{\,m+1}$, and together with the reverse inclusion we obtain
\[
\mathcal{I}_{p}^{(m)}=\mathbf{m}_{p}^{\,m+1}.
\]
\end{proof}

This allows us to characterize jets at a given point through derivatives and to determine the dimension of the corresponding fibers.

\begin{corollary}\label{cor:jets-caracterizacion-isomorfismo}
Let $M$ be a smooth manifold with or without boundary, let $p\in M$, and let $m\ge 0$. Fix a local chart centered at $p$ and write $D^\alpha$ for derivatives in those coordinates. Then, for any $f,g\in C^{\infty}(M)$ we have
\[
j^{m}_{p}f=j^{m}_{p}g
\quad\Longleftrightarrow\quad
D^{\alpha}f(p)=D^{\alpha}g(p),
\qquad\forall\,|\alpha|\le m.
\]
In particular, the chosen chart determines an isomorphism between the jet space $J^{m}_{p}(M)$ and the vector space
\[
\mathbb{K}^{N},
\qquad
N:=\binom{n+m}{m},
\]
through the map
\[
J^{m}_{p}(M)\longrightarrow \mathbb{K}^{N},
\qquad
j^{m}_{p}f\longmapsto \bigl(D^{\alpha}f(p)\bigr)_{|\alpha|\le m}.
\]
Under this identification, $\dim J^{m}_{p}(M)=\binom{n+m}{m}$.
\end{corollary}

\begin{proof}
By definition, $j^{m}_{p}f=j^{m}_{p}g$ if and only if the germs $[f]$ and $[g]$ satisfy
\[
[f]-[g]\in \mathbf{m}_{p}^{\,m+1}.
\]
By Lemma~\ref{lema:igualdad-ideales-jets}, $\mathbf{m}_{p}^{\,m+1}=\mathcal{I}_{p}^{(m)}$, that is,
\[
[f]-[g]\in \mathcal{I}_{p}^{(m)}
\quad\Longleftrightarrow\quad
D^{\alpha}(f-g)(p)=0,\ \forall\,|\alpha|\le m.
\]
This is equivalent to
\[
D^{\alpha}f(p)=D^{\alpha}g(p),
\qquad\forall\,|\alpha|\le m,
\]
which proves the first assertion.

To establish the isomorphism, use the smooth chart $(U,\phi)$ fixed in the statement, with $\phi(p)=0$. Define
\[
\Psi\colon J_{p}^{m}(M)\longrightarrow \mathbb{K}^{N},\qquad \Psi\left(j^{m}_{p}f\right)=\bigl(D^{\alpha}f(p)\bigr)_{|\alpha|\le m}
\]
This is well defined by the preceding criterion: if $j^{m}_{p}f=j^{m}_{p}g$, then $D^{\alpha}f(p)=D^{\alpha}g(p)$ for every $|\alpha|\le m$. It is also linear by linearity of differentiation.

To show that $\Psi$ is injective, Lemma~\ref{lema:igualdad-ideales-jets} gives:
\[
j^{m}_{p}f=j^{m}_{p}0
\quad\Longleftrightarrow\quad
[f]\in \mathbf{m}_{p}^{\,m+1} \quad\Longleftrightarrow\quad
D^{\alpha}f(p)=0,\text{ }\forall |\alpha|\leq m\quad\Longleftrightarrow \Psi\left(j_{p}^{m}f\right)=0.
\]

To show that $\Psi$ is surjective, given a collection of numbers \(\{a_{\alpha}\}_{|\alpha|\le m}\), consider the polynomial
\[
P(x)=\displaystyle\sum_{|\alpha|\le m}\frac{a_{\alpha}}{\alpha!}\,x^{\alpha},
\]
which satisfies $D^{\alpha}P(0)=a_{\alpha}$ for every multi-index $|\alpha|\leq m$. Define $f:=P\circ \phi\colon U\longrightarrow \mathbb{K}$. Since $p\in U$ and $f\in C^{\infty}(U)$, Theorem~\ref{extension de funciones suaves} gives an open set $W\subseteq M$ such that $p\in W\subseteq\overline{W}\subseteq U$ and a function $\widetilde{f}\in C^{\infty}(M)$ such that $\widetilde{f}\restriction_{W}=f$. Since $D^{\alpha}f(p)=D^{\alpha}(f\circ \phi^{-1})|_{\phi(p)}=D^{\alpha}P(0)$ and $D^{\alpha}\widetilde{f}(p)=D^{\alpha}f(p)$, we then have \[
\bigl(D^{\alpha}\widetilde{f}(p)\bigr)_{|\alpha|\le m}=(a_{\alpha})_{|\alpha|\le m}.
\] with $\widetilde{f}\in C^{\infty}(M)$.

Therefore, $\Psi$ is a linear isomorphism and the dimension of $J^{m}_{p}(M)$ is the number of multi-indices $\alpha\in\mathbb{N}_{0}^{n}$ with $|\alpha|\le m$, namely,
\[
\dim J^{m}_{p}(M)=\binom{n+m}{m}.
\]
\end{proof}
To gain more practice with differential operators between vector bundles, we extend the concept of a partial derivative to local sections of smooth vector bundles.

\begin{lemma}\label{lema: conmutador derivadas}
Let $\pi_{\mathbf{E}}\colon \mathbf{E}\longrightarrow M$ be a smooth vector bundle of rank $r$ over a smooth manifold $M$ with or without boundary, and let $(U,\phi)$ be a chart with local coordinates $\phi=(x^{1},\dots,x^{n})$ such that $U$ is the domain of a local trivialization
\[
\tau_{\mathbf{E}}\colon \mathbf{E}\restriction_{U}\longrightarrow U\times\mathbb{K}^{r}.
\]

For each multi-index $\beta\in\mathbb{N}_{0}^{n}$ define $D^{\beta}\colon \Gamma(\mathbf{E}\restriction_{U})\longrightarrow\Gamma(\mathbf{E}\restriction_{U})$ by
\[
(D^{\beta}\mathbf{u})(p)
:=
\tau_{\mathbf{E}}^{-1}\!\left(p,\,
D^{\beta}\bigl((\pi_{2}\circ\tau_{\mathbf{E}}\circ \mathbf{u})\bigr)(p)\right),
\qquad
\mathbf{u}\in\Gamma(\mathbf{E}\restriction_{U}),\;\; p\in U,
\]
where $D^{\beta}$ on the right-hand side is the partial derivative of order $|\beta|$ applied componentwise to the function \[
\pi_{2}\circ\tau_{\mathbf{E}}\circ \mathbf{u}\colon U\longrightarrow\mathbb{K}^{r}.
\]

Then, for every multi-index $\beta$ with $|\beta|\ge1$, we have
\[
[D^{\beta},f](\mathbf{u})
=
\tau_{\mathbf{E}}^{-1}\!\left(p,\,
\displaystyle\sum_{0<\gamma\le\beta}
\binom{\beta}{\gamma}
(D^{\gamma}f)(p)\,
D^{\beta-\gamma}\bigl((\pi_{2}\circ\tau_{\mathbf{E}}\circ \mathbf{u})\bigr)(p)
\right),
\]
for every $f\in C^{\infty}(U)$ and every $\mathbf{u}\in\Gamma(\mathbf{E}\restriction_{U})$.
\end{lemma}

\begin{proof}
Let $\mathbf{u}\in\Gamma(\mathbf{E}\restriction_{U})$. Its local representation is given by the function
\[
\widetilde u := \pi_{2}\circ\tau_{\mathbf{E}}\circ \mathbf{u}\colon U\longrightarrow\mathbb{K}^{r},
\]
so that
\[
\tau_{\mathbf{E}}(\mathbf{u}(p))=(p,\widetilde u(p)),
\qquad p\in U.
\]

In these coordinates, $D^{\beta}$ acts componentwise on $\widetilde u$.

\medskip

The commutator is defined by
\[
[D^{\beta},f]\mathbf{u}
=
D^{\beta}(f\mathbf{u})\;-\;f\,D^{\beta}\mathbf{u}.
\]
Applying $\pi_{2}\circ\tau_{\mathbf{E}}$ gives
\begin{equation}\label{eq:conmutador-expl}
(\pi_{2}\circ\tau_{\mathbf{E}}\circ[D^{\beta},f]\mathbf{u})(p)
=
D^{\beta}\!\bigl(f\widetilde u\bigr)(p)
-
f(p)\,D^{\beta}\widetilde u(p).
\end{equation}

\medskip

By the Leibniz rule (Proposition~\ref{prop: propiedades parciales multiindices incluye leibniz multiindices}), for every $p\in U$,
\begin{equation}\label{eq:leibniz-multinomial}
D^{\beta}(f\widetilde u)(p)
=
\displaystyle\sum_{\gamma\le\beta}
\binom{\beta}{\gamma}(D^{\gamma}f)(p)\;
D^{\,\beta-\gamma}\widetilde u(p).
\end{equation}
The term with $\gamma=0$ is
\[
(D^{0}f)(p)\,D^{\beta}\widetilde u(p)
=
f(p)\,D^{\beta}\widetilde u(p).
\]

Substituting \eqref{eq:leibniz-multinomial} into \eqref{eq:conmutador-expl} and subtracting the term $\gamma=0$ gives
\begin{equation}\label{eq:representacion-final}
(\pi_{2}\circ\tau_{\mathbf{E}}\circ[D^{\beta},f]\mathbf{u})(p)
=
\displaystyle\sum_{0<\gamma\le\beta}
\binom{\beta}{\gamma}
(D^{\gamma}f)(p)\,
D^{\,\beta-\gamma}\widetilde u(p).
\end{equation}

\medskip

Finally, to recover the bundle section from its representation in $\mathbb{K}^{r}$, apply the inverse of the trivialization:
\[
[D^{\beta},f]\mathbf{u}(p)
=
\tau_{\mathbf{E}}^{-1}\!\left(
p,\,
(\pi_{2}\circ\tau_{\mathbf{E}}\circ[D^{\beta},f]\mathbf{u})(p)
\right).
\]
Substituting \eqref{eq:representacion-final} gives
\[
[D^{\beta},f]\mathbf{u}
=
\tau_{\mathbf{E}}^{-1}\!\left(p,\,
\displaystyle\sum_{0<\gamma\le\beta}
\binom{\beta}{\gamma}
(D^{\gamma}f)(p)\,
D^{\,\beta-\gamma}\widetilde u(p)
\right),
\]
which is the desired formula.
\end{proof}

We continue with a lemma that will help characterize PDOs locally using what we have already established about jets.

\begin{lemma}[Vanishing on germs of order $\ge m{+}1$]\label{lema:anulacion-m+1}\index{vanishing on germs of given order}
Let $M$ be a smooth manifold with or without boundary and let $p\in M$. Consider the trivial bundle
\[
\underline{\mathbb{K}}^{r}=M\times\mathbb{K}^{r},
\qquad
\underline{\mathbb{K}}^{s}=M\times\mathbb{K}^{s},
\]
and let
\[
P\in \mathbf{PDO}^{(m)}\bigl(\underline{\mathbb{K}}^{r},
\underline{\mathbb{K}}^{s}\bigr).
\]

Define the map
\[
L_{p}\colon \ \bigl(C_{p}^{\infty}(M)\bigr)^{r}\longrightarrow \mathbb{K}^{s},
\qquad
L_{p}([u]) := (P(u))(p),
\]
which is well defined by locality (Remark~\ref{remark: pdo es local implica igualdad de pdos si hay igualdad de funciones}).

We write
\[
\mathbf{m}_{p}^{\,m+1}\bigl(C_{p}^{\infty}(M)\bigr)^{r}
:=
\left\{
\displaystyle\sum_{j=1}^{N}[f_{j}]\,[u_{j}]
\ \Bigg|\
[f_{j}]\in\mathbf{m}_{p}^{\,m+1},\
[u_{j}]\in\bigl(C_{p}^{\infty}(M)\bigr)^{r},\
N\in\mathbb{N}
\right\},
\]
 Here we recall that $(C_{p}^{\infty}(M))^{r}$ has the structure of a module over the ring $C_{p}^{\infty}(M)$, and note that $\mathbf{m}_{p}^{\,m+1}\bigl(C_{p}^{\infty}(M)\bigr)^{r}$ is an $C_{p}^{\infty}(M)$--submodule of $(C_{p}^{\infty}(M))^{r}$, so the quotient \[\
\bigl(C_{p}^{\infty}(M)\bigr)^{r}
\big/
\mathbf{m}_{p}^{\,m+1}\bigl(C_{p}^{\infty}(M)\bigr)^{r}\] has the structure of an $C_{p}^{\infty}(M)$--module. Since the field $\mathbb{K}$ can be regarded as a subset of $C_{p}^{\infty}(M)$, this space also inherits an $\mathbb{K}$--vector space structure.

With this notation and these hypotheses,
\[
L_{p}\bigl(
\mathbf{m}_{p}^{\,m+1}(C_{p}^{\infty}(M))^{r}
\bigr)
=
\{0\}
\subset\mathbb{K}^{s},
\]
and therefore $L_{p}$ induces a well-defined $\mathbb{K}$--linear map
\[
\overline{L}_{p}:\
\bigl(C_{p}^{\infty}(M)\bigr)^{r}
\big/
\mathbf{m}_{p}^{\,m+1}\bigl(C_{p}^{\infty}(M)\bigr)^{r}
\longrightarrow
\mathbb{K}^{s}.
\]
\end{lemma}

\begin{proof}
Let $[u]\in\bigl(C_{p}^{\infty}(M)\bigr)^{r}$. By $\mathbb{K}$--linearity of $P$, it suffices to prove that
\[
L_{p}\bigl([f]\,[u]\bigr)=0
\qquad\text{for every }[f]\in\mathbf{m}_{p}^{\,m+1}.
\]

By Lemma~\ref{lema:igualdad-ideales-jets},
\[
\mathbf{m}_{p}^{\,m+1}
=
\mathcal{I}_{p}^{(m)}
=
\bigl\{
[f]\in C_{p}^{\infty}(M)
\ \mid\
D^{\alpha}f(p)=0,\ |\alpha|\le m
\bigr\}.
\]

Take a smooth chart $(U,\phi)$ with $p\in U$, $\phi(U)$ convex, $\phi(p)=0$, and coordinate functions $x^{1},\dots,x^{n}$. Let $\tilde{f}:=f\circ\phi^{-1}$. The conditions $D^{\alpha}f(p)=0$ are equivalent to $D^{\alpha}\tilde{f}(0)=0$ for every $|\alpha|\le m$.

Taylor's theorem gives the identity
\[
\tilde{f}(x)
=
\displaystyle\sum_{|\alpha|=m+1}x^{\alpha}\,\tilde{g}_{\alpha}(x),
\qquad
\tilde{g}_{\alpha}\in C^{\infty}(\phi(U)).
\]
Let
\[
g_{\alpha}:=\tilde{g}_{\alpha}\circ\phi\in C^{\infty}(U),
\qquad |\alpha|=m+1.
\]
Then on $U$
\[
f=\displaystyle\sum_{|\alpha|=m+1}x^{\alpha}g_{\alpha}.
\]

We now extend all these functions simultaneously: by Theorem~\ref{extension de funciones suaves} there exist an open set
\[
p\in W\subseteq\overline{W}\subseteq U
\]
and smooth functions defined on $M$
\[
\tilde{x}^{i},\ \tilde{g}_{\alpha}^{\,M}\in C^{\infty}(M)
\]
such that
\[
\tilde{x}^{i}\restriction_{W}=x^{i}\restriction_{W},\qquad
\tilde{g}_{\alpha}^{\,M}\restriction_{W}=g_{\alpha}\restriction_{W}.
\]

Therefore, on $W$,
\[
f = \displaystyle\sum_{|\alpha|=m+1} x^{\alpha} g_{\alpha}
 = \displaystyle\sum_{|\alpha|=m+1} \tilde{x}^{\alpha}\,\tilde{g}_{\alpha}^{\,M},
\]
and passing to germs at $p$ gives
\[
[f]=\displaystyle\sum_{|\alpha|=m+1}[\tilde{x}^{\alpha}]\,[\tilde{g}_{\alpha}^{\,M}]
\qquad\text{in }C_{p}^{\infty}(M).
\]

Now let
\[
v_{\alpha} := \tilde{g}_{\alpha}^{\,M} u\in (C^{\infty}(M))^{r}.
\]

By linearity, it suffices to prove that
\[
L_{p}\bigl([\tilde{x}^{\alpha}][v_{\alpha}]\bigr)=0.
\]
We write
\[
\tilde{x}^{\alpha}
=
(\tilde{x}^{1})^{\alpha_{1}}\cdots(\tilde{x}^{n})^{\alpha_{n}}.
\]

Define inductively, for $k\ge1$,
\[
\operatorname{ad}(f)^{k}(P)
:=
\operatorname{ad}(f)\bigl(\operatorname{ad}(f)^{k-1}(P)\bigr),
\qquad
\operatorname{ad}(f)^{0}(P):=P.
\]

Apply the identity
\[
P(fw)=f\,P(w)+[P,f](w),
\qquad f\in C^{\infty}(M),
\]
successively $|\alpha|=m+1$ times, always using the global functions $\tilde{x}^{i}$. Upon evaluation at $p$, all terms retaining any factor $\tilde{x}^{i}$ vanish because $\tilde{x}^{i}(p)=0$, and only the term accumulating a commutator at each step survives.

Thus we obtain
\[
P\bigl(\tilde{x}^{\alpha}v_{\alpha}\bigr)(p)
=
\Bigl(
\operatorname{ad}(\tilde{x}^{1})^{\alpha_{1}}\cdots
\operatorname{ad}(\tilde{x}^{n})^{\alpha_{n}}P
\Bigr)\bigl(v_{\alpha}\bigr)(p).
\]

But this is an iterated commutator of length
\[
\alpha_{1}+\cdots+\alpha_{n}=|\alpha|=m+1,
\]
and vanishes by Proposition~\ref{equivalencia pdo con ad}. In particular,
\[
L_{p}\bigl([\tilde{x}^{\alpha}][v_{\alpha}]\bigr)=P\bigl(\tilde{x}^{\alpha}v_{\alpha}\bigr)(p)=0.
\]

Therefore,
\[
L_{p}([f][u])=0
\qquad\forall\,[f]\in\mathbf{m}_{p}^{\,m+1}.
\]

Now let
\[
L_{p}\colon \bigl(C_{p}^{\infty}(M)\bigr)^{r}\longrightarrow \mathbb{K}^{s}.
\]
Since
\[
L_{p}\bigl(\mathbf{m}_{p}^{\,m+1}(C_{p}^{\infty}(M))^{r}\bigr)=\{0\},
\]
if $[u]$ and $[w]$ represent the same class in the quotient, then $[u]-[w]\in\mathbf{m}_{p}^{\,m+1}(C_{p}^{\infty}(M))^{r}$ and therefore
\[
L_{p}([u])-L_{p}([w])=L_{p}([u]-[w])=0.
\]

We therefore define
\[
\overline{L}_{p}([u]):=L_{p}([u])=(P(u))(p).
\]
The definition is independent of the representative, and the $\mathbb{K}$--linearity of $L_{p}$ shows that $\overline{L}_{p}$ is $\mathbb{K}$--linear.

\end{proof}

\begin{remark}\label{obs:operadores-diferenciales-en-haces-lema-lema-anulacion-m-1-expresa}
Lemma~\ref{lema:anulacion-m+1} expresses the fact that a differential operator of order $\leq m$ between trivial bundles depends only on the $m$-jet of the germ evaluated at the point. Consequently, germs of order $\ge m{+}1$ are invisible to any differential operator of order $\leq m$ at the pointwise level, as the following example shows.
\end{remark}
\begin{example}\label{ej:operadores-diferenciales-en-haces-operador-modo-obtiene-evaluar-tanto-cual}
Let $M=\mathbb{R}$, $p=0$, and consider the operator
\[
P=\frac{d^{2}}{dx^{2}}\in\mathbf{PDO}^{(2)}\bigl(\underline{\mathbb{R}},\underline{\mathbb{R}}\bigr).
\]
Let $f(x)=x^{3}$, so that
\[
f(0)=f'(0)=f''(0)=0,
\text{ so that }[f]\in\mathbf{m}_{0}^{3}.
\]
Let $u(x)=\sin x$. Then
\[
P(fu)=\frac{d^{2}}{dx^{2}}(x^{3}\sin x),
\]
and we obtain
\[
(x^{3}\sin x)''
=
6x\sin x
+
6x^{2}\cos x
-
x^{3}\sin x.
\]
Evaluating at $p=0$ gives
\[
P(fu)(0)=0.
\]

Therefore,
\[
L_{0}\bigl([f][u]\bigr)=(P(fu))(0)=0,
\]
which illustrates Lemma~\ref{lema:anulacion-m+1}: every germ vanishing to order $m{+}1$ is invisible to a differential operator of order $\leq m$ when evaluated at the point.
\end{example}

\begin{theorem}[Local characterization of partial differential operators in the trivial case]\label{thm:pdo-expresion-local-trivial}\index{partial differential operator!local characterization!trivial case}
Let $M$ be a smooth manifold with or without boundary and consider the trivial bundles
\[
\mathbf{E}=\underline{\mathbb{K}}^{r}=M\times\mathbb{K}^{r},
\qquad
\mathbf{F}=\underline{\mathbb{K}}^{s}=M\times\mathbb{K}^{s}.
\]
Let
\[
P\colon C^{\infty}(M,\mathbb{K}^{r})\longrightarrow C^{\infty}(M,\mathbb{K}^{s})
\]
be an $\mathbb{K}$--linear map. Fix $m\ge 0$. The following are equivalent:

\begin{enumerate}[label=(\alph*)]

\item
$P\in\mathbf{PDO}^{(m)}(\mathbf{E},\mathbf{F})$.

\item
For every smooth chart $(U,\phi)$ with local coordinates $\phi=(x^{1},\dots,x^{n})$ there exist smooth functions
\[
a_{\alpha}\colon U\longrightarrow \operatorname{Hom}(\mathbb{K}^{r},\mathbb{K}^{s}),
\qquad |\alpha|\le m,
\]
such that, for every function $u\in C^{\infty}(M,\mathbb{K}^{r})$ and every $p\in U$,
\[
(Pu)(p)
=
\displaystyle\sum_{|\alpha|\le m}
a_{\alpha}(p)\,D^{\alpha}u(p),
\]
where $D^{\alpha}u(p)\in\mathbb{K}^{r}$ denotes the vector of partial derivatives of order $|\alpha|$ of the components of $u$ in the chart $(U,\phi)$.
\end{enumerate}
\end{theorem}

\begin{proof}
\noindent$(b)\Rightarrow(a)$:
Suppose that (b) holds. This means that, for each smooth chart $(U,\phi)$, there is a collection of smooth coefficients $a_{\alpha}$ such that, for every function $u\in C^{\infty}(M,\mathbb{K}^{r})$, the restriction of $Pu$ to $U$ satisfies
\[
(Pu)\restriction_{U} = \displaystyle\sum_{|\alpha|\le m} a_{\alpha}\,D^{\alpha}(u\restriction_{U}).
\]
Denote by $P_{U}$ the differential operator defined on the open set $U$ by this local expression:
\[
P_{U}\colon C^{\infty}(U,\mathbb{K}^{r}) \longrightarrow C^{\infty}(U,\mathbb{K}^{s}),
\qquad
v \longmapsto \displaystyle\sum_{|\alpha|\le m} a_{\alpha}\,D^{\alpha}v.
\]
We know that, on the domain $U$, each term $D^{\alpha}$ is a differential operator of order $|\alpha|$ (Corollary~\ref{derivada parcial alpha esima es pdo}), and multiplication by $a_{\alpha}$ has order $0$. By Proposition~\ref{prop: composicion de pdos es pdo}, we conclude that $P_{U}$ is a partial differential operator of order $\le m$ on $U$; that is,
\[
P_{U} \in \mathbf{PDO}^{(m)}(\mathbf{E}\restriction_{U},\mathbf{F}\restriction_{U}).
\]

To prove that the global operator $P$ belongs to $\mathbf{PDO}^{(m)}(\mathbf{E},\mathbf{F})$, we use the recursive definition based on commutators and adjoints. We must check that for any $m+1$ global smooth functions $f_{0}, \dots, f_{m} \in C^{\infty}(M)$, the iterated commutator vanishes:
\[
\operatorname{ad}(f_{0}) \cdots \operatorname{ad}(f_{m})(P) = 0.
\]
Consider $u \in C^{\infty}(M,\mathbb{K}^{r})$. Since we explicitly know the behavior of $Pu$ on $U$, for any $f \in C^{\infty}(M)$ we have the identity:
\[
(\operatorname{ad}(f)(P)u)\restriction_{U}
= ([P, f]u)\restriction_{U}
= P_{U}((fu)\restriction_{U}) - f\restriction_{U} P_{U}(u\restriction_{U})
= (\operatorname{ad}(f\restriction_{U})(P_{U}))(u\restriction_{U}).
\]
Applying this recursively to the iterated commutators gives, in each chart $(U,\phi)$:
\[
\left(\bigl(\operatorname{ad}(f_{0}) \cdots \operatorname{ad}(f_{m})(P)\bigr)(u)\right)\restriction_{U}
=
\bigl(\operatorname{ad}(f_{0}\restriction_{U}) \cdots \operatorname{ad}(f_{m}\restriction_{U})(P_{U})\bigr)(u\restriction_{U}).
\]
Since $P_{U}$ is a differential operator of order $\le m$ on $U$, its commutator with $m+1$ functions vanishes identically on $U$. Therefore,
\[
\left(\bigl(\operatorname{ad}(f_{0}) \cdots \operatorname{ad}(f_{m})(P)\bigr)(u)\right)\restriction_{U} = 0.
\]
Since we can choose an open cover of $M$ consisting of chart domains, and the function \[\left(\bigl(\operatorname{ad}(f_{0}) \cdots \operatorname{ad}(f_{m})(P)\bigr)(u)\right)\restriction_{U}\] vanishes on each of those domains,
\[
\bigl(\operatorname{ad}(f_{0}) \cdots \operatorname{ad}(f_{m})(P)\bigr)(u) = 0 \quad \text{throughout } M.
\]
Since $u$ was arbitrary, the operator $\operatorname{ad}(f_{0}) \dots \operatorname{ad}(f_{m})(P)$ is the zero operator. This proves that $P \in \mathbf{PDO}^{(m)}(\mathbf{E},\mathbf{F})$.

\medskip

$(a)\Rightarrow(b)$:
Now suppose that $P\in\mathbf{PDO}^{(m)}(\mathbf{E},\mathbf{F})$. In particular, by Remark~\ref{remark: pdo es local implica igualdad de pdos si hay igualdad de funciones}, $P$ is local. Since $\mathbf{E}$ and $\mathbf{F}$ are trivial bundles, a section of $\mathbf{E}$ can be viewed as a function $u\colon M\longrightarrow\mathbb{K}^{r}$ and a section of $\mathbf{F}$ as a function $v\colon M\longrightarrow\mathbb{K}^{s}$. For each $p\in M$ consider the map
\[
L_{p}\colon \ (C_{p}^{\infty}(M))^{r}\longrightarrow \mathbb{K}^{s},
\qquad
L_{p}([u]) := (Pu)(p),
\]
where $[u]$ denotes the germ of $u$ at $p$. Locality of $P$ ensures that $L_{p}$ is well defined.

By Lemma~\ref{lema:anulacion-m+1}, applied to $M$, the point $p$, and the operator $P$, we have
\[
L_{p}\bigl(\mathbf{m}_{p}^{\,m+1}(C_{p}^{\infty}(M))^{r}\bigr)=\{0\}.
\]

Therefore, $L_{p}$ induces a well-defined $\mathbb{K}$--linear map
\[
\overline{L}_{p}:\
(C_{p}^{\infty}(M))^{r}\big/
\mathbf{m}_{p}^{\,m+1}(C_{p}^{\infty}(M))^{r}
\longrightarrow \mathbb{K}^{s}.
\]

\medskip

We now justify in some detail the structure of the quotient
\[
(C_{p}^{\infty}(M))^{r}\big/
\mathbf{m}_{p}^{\,m+1}(C_{p}^{\infty}(M))^{r}.
\]
By definition of the submodule
\[
\mathbf{m}_{p}^{\,m+1}(C_{p}^{\infty}(M))^{r}
=
\left\{
\displaystyle\sum_{j=1}^{N}[f_{j}]\,[u_{j}]
\ \Bigg|\
[f_{j}]\in\mathbf{m}_{p}^{\,m+1},\
[u_{j}]\in(C_{p}^{\infty}(M))^{r},\text{ }N\in \mathbb{N}
\right\},
\]
two germs
\([u],[v]\in (C_{p}^{\infty}(M))^{r}\), with coordinates
\[
u=(u^{1},\dots,u^{r}),\qquad v=(v^{1},\dots,v^{r}),
\]
represent the same class in the quotient if and only if their difference belongs to that submodule, that is, if and only if
\[
[u^{j}]-[v^{j}]\in\mathbf{m}_{p}^{\,m+1}
\qquad\text{for each }j\in\{1,\dots,r\}.
\]

In other words,
\[
[u]=[v]\ \text{in}\
(C_{p}^{\infty}(M))^{r}\big/
\mathbf{m}_{p}^{\,m+1}(C_{p}^{\infty}(M))^{r}
\quad\Longleftrightarrow\quad
[u^{j}]=[v^{j}]\ \text{in}\
C_{p}^{\infty}(M)/\mathbf{m}_{p}^{\,m+1},
\ \forall j\in \{1,\dots,r\}.
\]

It follows that the quotient, viewed as a vector space,
\[
(C_{p}^{\infty}(M))^{r}\big/
\mathbf{m}_{p}^{\,m+1}(C_{p}^{\infty}(M))^{r}
\]
is isomorphic to the product of vector spaces
\[
\bigl(C_{p}^{\infty}(M)/\mathbf{m}_{p}^{\,m+1}\bigr)^{r}.
\]

Now take a smooth chart $(U,\phi)$ with $p\in U$ and coordinates $\phi=(x^{1},\dots,x^{n})$. By Definition~\ref{def:jet-funciones}, the space of jets of order $m$ at $p$ is
\[
J_{p}^{m}(M): =C_{p}^{\infty}(M)\big/\mathbf{m}_{p}^{\,m+1},
\]
and Corollary~\ref{cor:jets-caracterizacion-isomorfismo} states that each class $j^{m}_{p}f$ is determined exactly by the derivatives
\[
\bigl(D^{\alpha}f(p)\bigr)_{|\alpha|\le m},
\]
in coordinates $x^{1},\dots,x^{n}$, and that the map
\[
J^{m}_{p}(M)\longrightarrow \mathbb{K}^{N},
\qquad
j^{m}_{p}f\longmapsto \bigl(D^{\alpha}f(p)\bigr)_{|\alpha|\le m},
\]
is an isomorphism of vector spaces, where $N=\binom{n+m}{m}$.

Combining these observations, we obtain a chain of isomorphisms

\[
(C_{p}^{\infty}(M))^{r}\big/
\mathbf{m}_{p}^{\,m+1}(C_{p}^{\infty}(M))^{r}
\;\cong\;
\bigl(C_{p}^{\infty}(M)/\mathbf{m}_{p}^{\,m+1}\bigr)^{r}
\;\cong\;
\bigl(J^{m}_{p}(M)\bigr)^{r}
\;\cong\;
\mathbb{K}^{\,r\binom{n+m}{m}},
\]
where the last identification sends the class of
\[
u=(u^{1},\dots,u^{r})
\]
to the finite collection of partial derivatives
\[
\bigl(D^{\alpha}u^{j}(p)\bigr)_{|\alpha|\le m,\;1\le j\le r}.
\]

In particular, under this isomorphism the functional $\overline{L}_{p}$ corresponds to a linear transformation
\[
\mathbb{K}^{\,r\binom{n+m}{m}}\longrightarrow\mathbb{K}^{s},
\]
uniquely determined by a collection of linear transformations
\[
a_{\alpha}(p)\in \operatorname{Hom}(\mathbb{K}^{r},\mathbb{K}^{s}),
\qquad |\alpha|\le m,
\]
such that, for every germ $[u]\in (C_{p}^{\infty}(M))^{r}$,
\[
\overline{L}_{p}([u])
=
\displaystyle\sum_{|\alpha|\le m} a_{\alpha}(p)\,D^{\alpha}u(p),
\]
where $D^{\alpha}u(p)\in\mathbb{K}^{r}$ is the vector of partial derivatives of order $|\alpha|$ of the components of $u$ in the chart $(U,\phi)$.

By definition of $\overline{L}_{p}$ and $L_{p}$, if $u\in C^{\infty}(M,\mathbb{K}^{r})$, the germ $[u]$ at $p$ satisfies
\[
(Pu)(p)=L_{p}([u])=\overline{L}_{p}([u])
=
\displaystyle\sum_{|\alpha|\le m} a_{\alpha}(p)\,D^{\alpha}u(p).
\]

This proves the local formula in (b) in the chart $(U,\phi)$:
\[
(Pu)(p)
=
\displaystyle\sum_{|\alpha|\le m}
a_{\alpha}(p)\,D^{\alpha}u(p),
\qquad p\in U.
\]

It remains to check that the coefficients $a_\alpha$ depend smoothly on the point. This property is local. Fix $p_0\in U$ and choose open sets $V\Subset W\Subset U$ with $p_0\in V$, together with a function $\chi\in C_c^\infty(U)$ equal to one on $W$. For each vector in the standard basis $e_j\in\mathbb K^r$ and each multi-index $\gamma$ with $|\gamma|\leq m$, extend the function
\[
u_{\gamma,j}(x):=\chi(x)(x^1)^{\gamma_1}\cdots(x^n)^{\gamma_n}e_j.
\]
globally. On $V$, the cutoff function contributes to none of the derivatives. The formula already obtained, applied first to $u_{0,j}$, shows that $a_0(\,\cdot\,)e_j=P(u_{0,j})$ on $V$ and therefore that $a_0$ is smooth there. Suppose by induction that $a_\beta$ is smooth on $V$ for every $|\beta|<\ell$. If $|\gamma|=\ell$, then
\[
P(u_{\gamma,j})
=\gamma!\,a_\gamma e_j
+\sum_{|\beta|<\ell}a_\beta D^\beta u_{\gamma,j}
\qquad\text{on }V.
\]
Indeed, the derivatives of order $\ell$ of the monomial $x^\gamma$ vanish except for the one corresponding to $\gamma$, which equals $\gamma!$. Thus,
\[
a_\gamma e_j
=\frac1{\gamma!}\left(
P(u_{\gamma,j})-\sum_{|\beta|<\ell}a_\beta D^\beta u_{\gamma,j}
\right)
\qquad\text{on }V,
\]
and the right-hand side is smooth by the induction hypothesis. Running through the standard basis and the orders $\ell=0,\dots,m$, we conclude that all the $a_\alpha$ are smooth near $p_0$. Since $p_0$ was arbitrary, they are smooth throughout $U$.
\end{proof}

Finally, we give the local characterization of partial differential operators on a smooth real or complex vector bundle.
\begin{theorem}[Local characterization of partial differential operators on vector bundles]\label{thm:pdo expresion local explicita}\index{partial differential operator!local characterization!on vector bundles}
Let $M$ be a smooth manifold with or without boundary. Let $\pi_{\mathbf{E}}\colon \mathbf{E}\longrightarrow M$ and $\pi_{\mathbf{F}}\colon \mathbf{F}\longrightarrow M$ be smooth vector bundles of ranks $r$ and $s$ (real or complex), and let $m\ge 0$. For an $\mathbb{K}$--linear map
\[
P\colon \Gamma(\mathbf{E})\longrightarrow \Gamma(\mathbf{F}),
\]
the following are equivalent:

\begin{enumerate}[label=(\alph*)]
\item $P\in \mathbf{PDO}^{(m)}(\mathbf{E},\mathbf{F})$.

\item For every smooth chart $(U,\phi)$ with local trivializations
\[
\tau_{\mathbf{E}}\colon \mathbf{E}\restriction_{U}\longrightarrow U\times\mathbb{K}^{r},
\qquad
\tau_{\mathbf{F}}\colon \mathbf{F}\restriction_{U}\longrightarrow U\times\mathbb{K}^{s},
\]
there exist smooth functions
\[
a_{\alpha}\colon U\longrightarrow \operatorname{Hom}(\mathbb{K}^{r},\mathbb{K}^{s}),
\qquad |\alpha|\le m,
\]
such that, for every section $\mathbf{u}\in\Gamma(\mathbf{E})$ with local representation in $U$ given by
\[
\widetilde u:=\pi_{2}\circ\tau_{\mathbf{E}}\circ (\mathbf{u}\restriction_{U})\in C^{\infty}(U,\mathbb{K}^{r}),
\]
we have, for each $p\in U$,
\[
P(\mathbf{u})(p)
=
\tau_{\mathbf{F}}^{-1}\!\left(
p,\
\displaystyle\sum_{|\alpha|\le m} a_{\alpha}(p)\,D^{\alpha}\widetilde u(p)
\right).
\]
\end{enumerate}
\end{theorem}

\begin{proof}
Fix a smooth chart $(U,\phi)$ with coordinates $\phi=(x^{1},\dots,x^{n})$ and local trivializations $\tau_{\mathbf{E}}$ and $\tau_{\mathbf{F}}$ defined on $U$. Given a global section $\mathbf{u} \in \Gamma(\mathbf{E})$, denote by $\widetilde u \in C^{\infty}(U, \mathbb{K}^r)$ the local representation of $\mathbf{u}$
\[
\widetilde u := \pi_{2}\circ\tau_{\mathbf{E}}\circ (\mathbf{u}\restriction_{U}).
\]
Given a global section $\mathbf{v}\in\Gamma(\mathbf{F})$, its representation in the trivialization of $\mathbf{F}$ is, by definition, $\widetilde v:=\pi_2\circ\tau_{\mathbf{F}}\circ(\mathbf{v}\restriction_U)$.

\medskip

\noindent$(b)\Rightarrow(a)$:
Suppose that (b) holds. Define an auxiliary operator $\widetilde{P}_{U}$ on the space of vector-valued functions on the open set $U$ by the explicit expression in the hypothesis:
\[
\widetilde{P}_{U}\colon C^{\infty}(U,\mathbb{K}^{r}) \longrightarrow C^{\infty}(U,\mathbb{K}^{s}),
\qquad
\widetilde{P}_{U}(w)=\displaystyle\sum_{|\alpha|\le m} a_{\alpha}\,D^{\alpha}w.
\]
By Theorem~\ref{thm:pdo-expresion-local-trivial}, the operator $\widetilde{P}_{U}$ is a partial differential operator of order $\leq m$ between the trivial bundles over $U$. Hypothesis (b) states that the global operator $P$ and this local operator $\widetilde{P}_{U}$ are related through the trivializations. Specifically, for every global section $\mathbf{u} \in \Gamma(\mathbf{E})$:

\[P(\mathbf{u})\restriction_{U}=\tau_{\mathbf{F}}^{-1}\left(
\ \cdot
 \ \ ,
\displaystyle\sum_{|\alpha|\le m} a_{\alpha}D^{\alpha}\widetilde u
\right),\]
so, composing with $\tau_{\mathbf{F}}$ and then projecting onto the second coordinate,
\[
\pi_{2}\circ\tau_{\mathbf{F}}\circ (P(\mathbf{u})\restriction_{U}) = \widetilde{P}_{U}(\widetilde u).
\]
To check that $P\in\mathbf{PDO}^{(m)}(\mathbf{E},\mathbf{F})$, it suffices to prove that
\[
\operatorname{ad}(f_{0})\cdots\operatorname{ad}(f_{m})(P)=0,
\qquad
\forall\, f_{0},\dots,f_{m}\in C^{\infty}(M).
\]
Let $\mathbf{u}\in\Gamma(\mathbf{E})$ and fix arbitrary functions $f_{0},\dots,f_{m}\in C^{\infty}(M)$. We work in the chart $(U,\phi)$ with the chosen trivializations.

\medskip

Given $f\in C^{\infty}(M)$ and $\mathbf{u}\in\Gamma(\mathbf{E})$, first observe that the pointwise action of $f$ passes to the local representation through the identity
\begin{equation}\label{eq:representacion-local-fu}
\widetilde{f\mathbf{u}}
=
\pi_{2}\circ\tau_{\mathbf{E}}\bigl((f\mathbf{u})\restriction_{U}\bigr)
=
(f\restriction_{U})\,\widetilde u.
\end{equation}
This equality follows directly from linearity of the trivialization $\tau_{\mathbf{E}}$ on each fiber.

\medskip

Now apply $\pi_{2}\circ\tau_{\mathbf{F}}$ to the commutator $\operatorname{ad}(f)(P)$:
\[
\pi_{2}\circ\tau_{\mathbf{F}}\circ\bigl(\operatorname{ad}(f)(P)(\mathbf{u})\bigr)\restriction_{U}
=
\pi_{2}\circ\tau_{\mathbf{F}}\Bigl( P(f\mathbf{u})\restriction_{U} - (fP(\mathbf{u}))\restriction_{U} \Bigr)
\]
\[
=
\widetilde{P}_{U}\bigl(\widetilde{f\mathbf{u}}\bigr)
-
(f\restriction_{U})\,\widetilde{P}_{U}(\widetilde u)
\]
\[
\overset{\eqref{eq:representacion-local-fu}}{=}
\widetilde{P}_{U}\bigl((f\restriction_{U})\,\widetilde u\bigr)
-
(f\restriction_{U})\,\widetilde{P}_{U}(\widetilde u)
\]
\[
=
\bigl[\widetilde{P}_{U},\, f\restriction_{U}\bigr](\widetilde u)
\]
\[
=
\bigl(\operatorname{ad}(f\restriction_{U})(\widetilde{P}_{U})\bigr)(\widetilde u).
\]

\medskip

Repeating this argument $m+1$ times and using at each step identity \eqref{eq:representacion-local-fu} together with the compatibility between $P$ and $\widetilde{P}_{U}$ provided by the trivializations gives
\[
\pi_{2}\circ\tau_{\mathbf{F}}\circ
\bigl(\operatorname{ad}(f_{0})\cdots\operatorname{ad}(f_{m})(P)(\mathbf{u})\bigr)
\restriction_{U}
=
\bigl(
\operatorname{ad}(f_{0}\restriction_{U})\cdots\operatorname{ad}(f_{m}\restriction_{U})
(\widetilde{P}_{U})
\bigr)(\widetilde u).
\]

\medskip

By Theorem~\ref{thm:pdo-expresion-local-trivial}, the operator $\widetilde{P}_{U}$ is a partial differential operator of order $\le m$ on the trivial bundle over $U$. By Proposition~\ref{equivalencia pdo con ad}, every iterated commutator of length $m+1$ with smooth functions on $U$ vanishes:
\[
\operatorname{ad}(f_{0}\restriction_{U})\cdots\operatorname{ad}(f_{m}\restriction_{U})
(\widetilde{P}_{U})
=
0.
\]
Therefore, the right-hand side of the preceding equality is $0$ for every $\widetilde u$, and consequently
\[
\operatorname{ad}(f_{0})\cdots\operatorname{ad}(f_{m})(P)(\mathbf{u})=0
\qquad\text{on }U.
\]

\medskip

Since $M$ admits an open cover by charts that are trivialization domains for $\mathbf{E}$ and $\mathbf{F}$, and the preceding identity holds on each of them, we conclude that
\[
\operatorname{ad}(f_{0})\cdots\operatorname{ad}(f_{m})(P)=0
\qquad\text{throughout }M.
\]
This proves that $P\in\mathbf{PDO}^{(m)}(\mathbf{E},\mathbf{F})$.

\medskip

\noindent$(a)\Rightarrow(b)$:

Now suppose that $P \in \mathbf{PDO}^{(m)}(\mathbf{E},\mathbf{F})$. First, $P$ is a local operator (Remark~\ref{remark: pdo es local implica igualdad de pdos si hay igualdad de funciones}). This means that if two global sections agree on an open set, their images under $P$ agree on that open set. This locality allows us to define the version of the operator in each trivialization.

We define an operator $\widetilde{P}\colon C^{\infty}(U, \mathbb{K}^r) \longrightarrow C^{\infty}(U, \mathbb{K}^s)$ transferring the structure of $P$ to $U$ through a pointwise construction. Let $w \in C^{\infty}(U, \mathbb{K}^r)$ and $p \in U$. To evaluate $\widetilde{P}(w)$ at $p$, we need to transfer $w$ to a global section on which $P$ can act. Choose an open set $V\subseteq M$ such that $p \in V$ and $\overline{V} \subseteq U$. By Proposition~\ref{funciones flan} (existence of cutoff functions), there exists a smooth function $\psi\colon M \longrightarrow \mathbb{R}$ such that $\psi \equiv 1$ on $\overline{V}$ and $\supp(\psi) \subseteq U$. Construct a global extension $\widehat{\mathbf{w}} \in \Gamma(\mathbf{E})$ defined by:
\begin{equation}\label{eq: extension suave de seccion para haz trivial}
 \widehat{\mathbf{w}}(q) :=
\begin{cases}
\tau_{\mathbf{E}}^{-1}(q, \psi(q)w(q)) & \text{if } q \in U, \\
0 & \text{if } q \in M \setminus U.
\end{cases}
\end{equation}
We then define the value of $\widetilde{P}(w)$ at the point $p$ by:
\[
\widetilde{P}(w)(p) := \pi_{2} \circ \tau_{\mathbf{F}} \bigl( P(\widehat{\mathbf{w}})(p) \bigr).
\]
By locality of $P$, this value is independent of the choice of cutoff function $\psi$ and open set $V$, since any other extension would agree with $\widehat{\mathbf{w}}$ on a neighborhood of $p$. Therefore, $\widetilde{P}$ is well defined as a linear operator between sections of the trivial bundles.

To apply Theorem~\ref{thm:pdo-expresion-local-trivial}, we must prove that $\widetilde{P}$ is a partial differential operator of order $\le m$ between the trivial bundles over $U$. That is, we must show that for any $g_0, \cdots, g_m \in C^\infty(U)$ we have:
\[
\operatorname{ad}(g_0)\cdots\operatorname{ad}(g_m)(\widetilde{P})=0.
\]
Let $p\in U$ be arbitrary. Choose an open set $V\subseteq U$ such that $\overline{V}\subseteq U$ and a smooth function $\psi\colon M\longrightarrow\mathbb{R}$ with $\psi\equiv 1$ on $\overline{V}$ and $\supp(\psi)\subseteq U$. This choice fixes a uniform procedure for extending scalar functions and sections defined on $U$.

\medskip

Given a function $g\in C^{\infty}(U)$, define its global extension $\widehat{g}\in C^{\infty}(M,\mathbb{K})$ by
\[
\widehat{g}(q)=
\begin{cases}
\psi(q)\,g(q), & q\in U,\\
0, & q\in M\setminus U.
\end{cases}
\]
By construction, $\widehat{g}\restriction_{V}=g\restriction_{V}$ and $\supp(\widehat{g})\subseteq U$. Likewise, given a vector-valued function $w\in C^{\infty}(U,\mathbb{K}^{r})$, define its global extension $\widehat{\mathbf{w}}\in\Gamma(\mathbf{E})$ as in \eqref{eq: extension suave de seccion para haz trivial}:
\[
\widehat{\mathbf{w}}(q)
=
\begin{cases}
\tau_{\mathbf{E}}^{-1}(q,\psi(q)\,w(q)), & q\in U,\\
0, & q\in M\setminus U.
\end{cases}
\]
In particular, $\widehat{\mathbf{w}}\restriction_{V}=\tau_{\mathbf{E}}^{-1}(\cdot,w(\cdot))$.

\medskip

With these extensions fixed, we relate the commutators of $\widetilde{P}$ to those of $P$. We begin with the case of a single commutator. Let $g\in C^{\infty}(U)$ and let $\widehat{g}\in C^{\infty}(M)$ be its global extension defined above. Let $w\in C^{\infty}(U,\mathbb{K}^{r})$, and let $\widehat{\mathbf{w}}\in\Gamma(\mathbf{E})$ be its global extension. By definition of $\widetilde{P}$,
\[
\widetilde{P}(w)(q)
=
\pi_{2}\circ\tau_{\mathbf{F}}\bigl(P(\widehat{\mathbf{w}})(q)\bigr),
\qquad q\in V.
\]

Now fix $q\in V$. Since $\psi\equiv 1$ on $V$, we have
\[
\widehat{g}(q)=g(q),
\qquad
\widehat{\mathbf{w}}(q)=\tau_{\mathbf{E}}^{-1}(q,w(q)).
\]
Let $\widehat{\mathbf{w}}_g$ be the global extension of the vector-valued function $gw$, constructed using the same function $\psi$:
\[
\widehat{\mathbf{w}}_g(q')
=
\begin{cases}
\tau_{\mathbf{E}}^{-1}\bigl(q',\psi(q')\,g(q')w(q')\bigr), & q'\in U,\\
0, & q'\notin U.
\end{cases}
\]
For $q'\in V$,
\[
\widehat{\mathbf{w}}_g(q')
=
\tau_{\mathbf{E}}^{-1}\bigl(q',g(q')w(q')\bigr)
=
\widehat{g}(q')\,\widehat{\mathbf{w}}(q'),
\]
where the last equality is scalar multiplication in the fiber.

Since $\widehat{\mathbf{w}}_g$ and $\widehat{g}\,\widehat{\mathbf{w}}$ agree on a neighborhood of $q$, locality of $P$ implies
\[
P(\widehat{\mathbf{w}}_g)(q)=P(\widehat{g}\,\widehat{\mathbf{w}})(q).
\]

Therefore,
\begin{align*}
\bigl(\operatorname{ad}(g)(\widetilde{P}) w\bigr)(q)
&=
\widetilde{P}(gw)(q) - g(q)\,\widetilde{P}(w)(q) \\
&=
\pi_{2}\circ\tau_{\mathbf{F}}\bigl(P(\widehat{\mathbf{w}}_g)(q)\bigr)
-
\widehat{g}(q)\,\pi_{2}\circ\tau_{\mathbf{F}}\bigl(P(\widehat{\mathbf{w}})(q)\bigr) \\
&=
\pi_{2}\circ\tau_{\mathbf{F}}\bigl(
P(\widehat{g}\,\widehat{\mathbf{w}})(q) - \widehat{g}(q)\,P(\widehat{\mathbf{w}})(q)
\bigr) \\
&=
\pi_{2}\circ\tau_{\mathbf{F}}\bigl(
\operatorname{ad}(\widehat{g})(P)(\widehat{\mathbf{w}})(q)
\bigr).
\end{align*}

Thus, for every $q\in V$,
\[
\bigl(\operatorname{ad}(g)(\widetilde{P}) w\bigr)(q)
=
\pi_{2}\circ\tau_{\mathbf{F}}\bigl(
\operatorname{ad}(\widehat{g})(P)(\widehat{\mathbf{w}})(q)
\bigr).
\]

Applying the same reasoning inductively gives, for any $g_{0},\dots,g_{m}\in C^{\infty}(U)$ and their extensions $\widehat{g}_{0},\dots,\widehat{g}_{m}$,
\[
\bigl(\operatorname{ad}(g_{0})\cdots\operatorname{ad}(g_{m})
(\widetilde{P}) w\bigr)(q)
=
\pi_{2}\circ\tau_{\mathbf{F}}\bigl(
\operatorname{ad}(\widehat{g}_{0})\cdots\operatorname{ad}(\widehat{g}_{m})(P)(\widehat{\mathbf{w}})
\bigr)(q),
\qquad q\in V.
\]

Since by hypothesis $P \in \mathbf{PDO}^{(m)}(\mathbf{E},\mathbf{F})$, $\operatorname{ad}(\hat{g}_{0})\cdots\operatorname{ad}(\widehat{g}_{m})(P)$ is the zero operator on $M$. In particular, $\bigl(
\operatorname{ad}(\widehat{g}_{0})\cdots\operatorname{ad}(\widehat{g}_{m})(P)(\widehat{\mathbf{w}})
\bigr)(q)=0$ and therefore $\bigl(\bigl(\operatorname{ad}(g_{0})\cdots\operatorname{ad}(g_{m})
(\widetilde{P})\bigr)( w)\bigr)(q)=0$. Taking $q=p$, $\bigl(\bigl(\operatorname{ad}(g_{0})\cdots\operatorname{ad}(g_{m})
(\widetilde{P})\bigr)( w)\bigr)(p)=0$, and using the fact that $p$ and $w$ were arbitrary, we conclude that
\[
\operatorname{ad}(g_0)\cdots\operatorname{ad}(g_m)(\widetilde{P}) = 0 \quad \text{throughout } U.
\]
Therefore, $\widetilde{P} \in \mathbf{PDO}^{(m)}(\underline{\mathbb{K}}^{r}, \underline{\mathbb{K}}^{s})$, where $\underline{\mathbb{K}}^{r}$ and $\underline{\mathbb{K}}^{s}$ are the trivial bundles of ranks $r$ and $s$, respectively, over $U$. Applying Theorem~\ref{thm:pdo-expresion-local-trivial} to the operator $\widetilde{P}$ ensures the existence of smooth functions $a_{\alpha}\colon U\longrightarrow\operatorname{Hom}(\mathbb{K}^{r},\mathbb{K}^{s})$ such that, for every $w \in C^\infty(U, \mathbb{K}^r)$ and every $p \in U$:
\[
\widetilde{P}(w)(p) = \displaystyle\sum_{|\alpha|\le m}a_{\alpha}(p)\,D^{\alpha}w(p).
\]
Finally, consider a global section $\mathbf{u} \in \Gamma(\mathbf{E})$. Let $\widetilde u \in C^\infty(U, \mathbb{K}^r)$ be its local representation. For any $p \in U$, $\mathbf{u}$ itself serves as a global extension of its local behavior. By definition of $\widetilde{P}$, we have for every $p\in U$:
\[
\pi_{2}\circ\tau_{\mathbf{F}}\bigl( P(\mathbf{u})(p) \bigr) = \widetilde{P}(\widetilde u)(p) = \displaystyle\sum_{|\alpha|\le m}a_{\alpha}(p)\,D^{\alpha}\widetilde u(p).
\]
Applying $\tau_{\mathbf{F}}^{-1}(p, \cdot)$ to both sides gives \[
P(\mathbf{u})(p)
=
\tau_{\mathbf{F}}^{-1}\!\left(
p,\
\displaystyle\sum_{|\alpha|\le m} a_{\alpha}(p)\,D^{\alpha}\widetilde u(p)
\right).
\]
\end{proof}

\begin{corollary}[Local expression of a PDO using local frames]\label{cor:pdo-expresion-con-marcos}\index{local expression of a PDO using local frames}
Let $M$ be a smooth manifold with or without boundary. Let $\pi_{\mathbf{E}}\colon \mathbf{E}\longrightarrow M$ and $\pi_{\mathbf{F}}\colon \mathbf{F}\longrightarrow M$ be smooth vector bundles of ranks $r$ and $s$. Under the hypotheses of Theorem~\ref{thm:pdo expresion local explicita}, for a continuous linear operator
\[
P\in\mathbf{Op}(\mathbf{E},\mathbf{F})
\]
and an integer $m\geq 0$, the following are equivalent:

\begin{enumerate}[label=(\alph*)]
\item $P\in\mathbf{PDO}^{(m)}(\mathbf{E},\mathbf{F})$.

\item For every smooth chart $(U,\phi)$ and every pair of smooth local frames
\[
\{\mathbf{e}_{1},\dots,\mathbf{e}_{r}\} \text{ in } \mathbf{E}\restriction_{U},
\qquad
\{\mathbf{f}_{1},\dots,\mathbf{f}_{s}\} \text{ in } \mathbf{F}\restriction_{U},
\]
there exist smooth functions
\[
(a_{\alpha})_{ij}\colon U\longrightarrow\mathbb{K},
\qquad
1\leq i\leq r,\ 1\leq j\leq s,\ |\alpha|\leq m,
\]
such that every section written in components as
\[
\mathbf{u}=\displaystyle\sum_{i=1}^{r}u^{i}\mathbf{e}_{i},
\qquad u^{i}\in C^{\infty}(U,\mathbb{K}),
\]
satisfies the identity
\[
P(\mathbf{u})
=
\displaystyle\sum_{j=1}^{s}
\left(
\displaystyle\sum_{|\alpha|\leq m}\displaystyle\sum_{i=1}^{r}
(a_{\alpha})_{ij}\,D^{\alpha}u^{i}
\right)
\mathbf{f}_{j}.
\]
on $U$.
\end{enumerate}
\end{corollary}

\begin{proof}
(a)$\Rightarrow$(b).
Fix a chart $(U,\phi)$ and local frames $\{\mathbf{e}_{i}\}$ in $\mathbf{E}\restriction_{U}$ and $\{\widehat{\mathbf{e}}_{j}\}$ in $\mathbf{F}\restriction_{U}$. Recall that, by Proposition~\ref{marcos asociados con trivializaciones}, a smooth local frame determines a smooth local trivialization
\[
\boldsymbol{\tau}_{\mathbf{E}}\colon \mathbf{E}\restriction_{U}\longrightarrow U\times\mathbb{K}^{r},
\qquad
\boldsymbol{\tau}_{\mathbf{F}}\colon \mathbf{F}\restriction_{U}\longrightarrow U\times\mathbb{K}^{s},
\]
characterized by
\[
\boldsymbol{\tau}_{\mathbf{E}}(\mathbf{e}_{i}(p))=(p,e_{i}^{0}),
\qquad
\boldsymbol{\tau}_{\mathbf{F}}(\widehat{\mathbf{e}}_{j}(p))=(p,\widehat{e}_{j}^{0}),
\]
where $\{e_{i}^{0}\}$ and $\{\widehat{e}_{j}^{0}\}$ are the standard bases of $\mathbb{K}^{r}$ and $\mathbb{K}^{s}$, respectively.

Applying Theorem~\ref{thm:pdo expresion local explicita} to these trivializations, there exist smooth maps
\[
a_{\alpha}\colon U\longrightarrow\operatorname{Hom}(\mathbb{K}^{r},\mathbb{K}^{s}),
\qquad |\alpha|\leq m,
\]
such that, for the local representation
\[
\widetilde u:=\pi_{2}\circ\boldsymbol{\tau}_{\mathbf{E}}\circ \mathbf{u}\in C^{\infty}(U,\mathbb{K}^{r}),
\]
we have
\[
\pi_{2}\circ\boldsymbol{\tau}_{\mathbf{F}}\bigl(P(\mathbf{u})\bigr)
=
\displaystyle\sum_{|\alpha|\leq m} a_{\alpha}\,D^{\alpha}\widetilde u.
\]

To obtain the component expression, write the action of $a_{\alpha}(p)$ on the standard basis:
\[
a_{\alpha}(p)(e_{i}^{0})
=
\displaystyle\sum_{j=1}^{s}(a_{\alpha})_{ij}(p)\,\widehat{e}_{j}^{0}.
\]
This defines smooth functions $(a_{\alpha})_{ij}$, which are precisely the matrix entries of $a_{\alpha}(p)$.

If $\mathbf{u}=\displaystyle\sum_{i=1}^{r}u^{i}\mathbf{e}_{i}$, then, for each $p\in U$,
\[
\mathbf{u}(p)=\displaystyle\sum_{i=1}^{r}u^{i}(p)\,\mathbf{e}_{i}(p).
\]
The trivialization induced by the local frame $\{\mathbf{e}_{1},\dots,\mathbf{e}_{r}\}$ is given by
\[
\boldsymbol{\tau}_{\mathbf{E}}(\mathbf{e}_{i}(p))=(p,e_{i}^{0}),
\]
where $(e_{i}^{0})$ is the standard basis of $\mathbb{K}^{r}$. By linearity on each fiber,
\[
\boldsymbol{\tau}_{\mathbf{E}}(\mathbf{u}(p))
=
\boldsymbol{\tau}_{\mathbf{E}}\!\left(\displaystyle\sum_{i=1}^{r}u^{i}(p)\mathbf{e}_{i}(p)\right)
=
\left(p,\displaystyle\sum_{i=1}^{r}u^{i}(p)e_{i}^{0}\right)
=
\bigl(p,(u^{1}(p),\dots,u^{r}(p))\bigr).
\]
Thus, the local representation
\[
\widetilde u=\pi_{2}\circ\boldsymbol{\tau}_{\mathbf{E}}\circ \mathbf{u}\colon U\longrightarrow\mathbb{K}^{r}
\]
is explicitly given by
\[
\widetilde u(p)=(u^{1}(p),\dots,u^{r}(p)).
\]

Therefore,
\[
a_{\alpha}(p)\,D^{\alpha}\widetilde u(p)
=
\left(
\displaystyle\sum_{i=1}^{r}(a_{\alpha})_{i1}(p)D^{\alpha}u^{i}(p),\dots,
\displaystyle\sum_{i=1}^{r}(a_{\alpha})_{is}(p)D^{\alpha}u^{i}(p)
\right).
\]
The local representation of $P(\mathbf{u})$ is
\[
\pi_{2}\circ\boldsymbol{\tau}_{\mathbf{F}}\bigl(P(\mathbf{u})\bigr)(p)
=
\left(
\displaystyle\sum_{|\alpha|\le m}\displaystyle\sum_{i=1}^{r}
(a_{\alpha})_{i1}(p)D^{\alpha}u^{i}(p),
\dots,
\displaystyle\sum_{|\alpha|\le m}\displaystyle\sum_{i=1}^{r}
(a_{\alpha})_{is}(p)D^{\alpha}u^{i}(p)
\right).
\]
Since $\boldsymbol{\tau}_{\mathbf{F}}(\mathbf{f}_{j}(p))=(p,f_{j}^{0})$, inverting $\boldsymbol{\tau}_{\mathbf{F}}$ gives
\[
P(\mathbf{u})(p)
=
\displaystyle\sum_{j=1}^{s}
\left(
\displaystyle\sum_{|\alpha|\le m}\displaystyle\sum_{i=1}^{r}
(a_{\alpha})_{ij}(p)D^{\alpha}u^{i}(p)
\right)
\widehat{\mathbf{e}}_{j}(p),
\]
and therefore on $U$
\[
P(\mathbf{u})
=
\displaystyle\sum_{j=1}^{s}
\left(
\displaystyle\sum_{|\alpha|\le m}\displaystyle\sum_{i=1}^{r}
(a_{\alpha})_{ij}D^{\alpha}u^{i}
\right)
\widehat{\mathbf{e}}_{j}.
\]

\medskip

(b)$\Rightarrow$(a).
Suppose that (b) holds. We must prove that $P\in\mathbf{PDO}^{(m)}(\mathbf{E},\mathbf{F})$ in the sense of Theorem~\ref{thm:pdo expresion local explicita}. Let $(U,\phi)$ be an arbitrary smooth chart of $M$ and take arbitrary smooth local trivializations
\[
\boldsymbol{\tau}_{\mathbf{E}}\colon \mathbf{E}\restriction_{U}\longrightarrow U\times\mathbb{K}^{r},
\qquad
\boldsymbol{\tau}_{\mathbf{F}}\colon \mathbf{F}\restriction_{U}\longrightarrow U\times\mathbb{K}^{s}.
\]

By Proposition~\ref{marcos asociados con trivializaciones}, each of these trivializations induces a smooth local frame over $U$. More precisely, there exists a smooth local frame $\{\mathbf{e}_{1},\dots,\mathbf{e}_{r}\}$ in $\mathbf{E}\restriction_{U}$ associated with $\boldsymbol{\tau}_{\mathbf{E}}$, characterized by
\[
\mathbf{e}_{i}(p)=\boldsymbol{\tau}_{\mathbf{E}}^{-1}(p,e_{i}^{0}),\qquad 1\leq i\leq r,
\]
where $\{e_{i}^{0}\}$ is the standard basis of $\mathbb{K}^{r}$. The trivialization $\boldsymbol{\tau}_{\mathbf{F}}$ induces the smooth local frame $\{\widehat{\mathbf{e}}_{1},\dots,\widehat{\mathbf{e}}_{s}\}$ in $\mathbf{F}\restriction_{U}$ given by
\[
\widehat{\mathbf{e}}_{j}(p)=\boldsymbol{\tau}_{\mathbf{F}}^{-1}(p,\widehat{e}_{j}^{0}),\qquad 1\leq j\leq s,
\]
where $\{\widehat{e}_{j}^{0}\}$ is the standard basis of $\mathbb{K}^{s}$.

Since by hypothesis (b) holds for every pair of local frames, in particular for the frames $\{\mathbf{e}_{i}\}$ and $\{\widehat{\mathbf{e}}_{j}\}$ associated with these trivializations, there exist smooth functions
\[
(a_{\alpha})_{ij}\colon U\longrightarrow\mathbb{K},
\qquad 1\leq i\leq r,\ 1\leq j\leq s,\ |\alpha|\leq m,
\]
such that, for every section
\[
\mathbf{u}=\displaystyle\sum_{i=1}^{r}u^{i}\mathbf{e}_{i},
\qquad u^{i}\in C^{\infty}(U,\mathbb{K}),
\]
the identity
\[
P(\mathbf{u})
=
\displaystyle\sum_{j=1}^{s}
\left(
\displaystyle\sum_{|\alpha|\leq m}\displaystyle\sum_{i=1}^{r}
(a_{\alpha})_{ij}\,D^{\alpha}u^{i}
\right)
\widehat{\mathbf{e}}_{j}.
\]
holds on $U$.

Now define, for each multi-index $\alpha$, a map
\[
a_{\alpha}\colon U\longrightarrow\operatorname{Hom}(\mathbb{K}^{r},\mathbb{K}^{s})
\]
as follows: for each $p\in U$ and each vector $v=(v^{1},\dots,v^{r})\in\mathbb{K}^{r}$, set
\[
a_{\alpha}(p)(v)
:=
\left(
\displaystyle\sum_{i=1}^{r}(a_{\alpha})_{i1}(p)v^{i},
\dots,
\displaystyle\sum_{i=1}^{r}(a_{\alpha})_{is}(p)v^{i}
\right)\in\mathbb{K}^{s}.
\]
In other words, the matrix of $a_{\alpha}(p)$ in the standard bases $\{e_{i}^{0}\}$ and $\{\widehat{e}_{j}^{0}\}$ is precisely
\[
\bigl((a_{\alpha})_{ij}(p)\bigr)_{1\leq i\leq r,\ 1\leq j\leq s}.
\]
Since the functions $(a_{\alpha})_{ij}$ are smooth, it follows that $a_{\alpha}$ is a smooth map with values in $\operatorname{Hom}(\mathbb{K}^{r},\mathbb{K}^{s})$.

We show that these $a_{\alpha}$ provide the local expression of $P$ in the trivializations $\boldsymbol{\tau}_{\mathbf{E}}$ and $\boldsymbol{\tau}_{\mathbf{F}}$. Let $\mathbf{u}\in\Gamma(\mathbf{E})$ and write on $U$
\[
\mathbf{u}=\displaystyle\sum_{i=1}^{r}u^{i}\mathbf{e}_{i},
\qquad u^{i}\in C^{\infty}(U,\mathbb{K}).
\]
By construction of $\boldsymbol{\tau}_{\mathbf{E}}$ from the frame $\{\mathbf{e}_{i}\}$, the local representation of $\mathbf{u}$ is
\[
\widetilde u:=\pi_{2}\circ\boldsymbol{\tau}_{\mathbf{E}}\circ \mathbf{u}\colon U\longrightarrow\mathbb{K}^{r},
\qquad
\widetilde u(p)=(u^{1}(p),\dots,u^{r}(p)).
\]
In particular,
\[
D^{\alpha}\widetilde u(p)
=
\bigl(D^{\alpha}u^{1}(p),\dots,D^{\alpha}u^{r}(p)\bigr).
\]

On the other hand, using the formula in (b) for $P(\mathbf{u})$ and applying $\boldsymbol{\tau}_{\mathbf{F}}$ gives
\[
\boldsymbol{\tau}_{\mathbf{F}}\bigl(P(\mathbf{u})(p)\bigr)
=
\left(
p,\
\displaystyle\sum_{j=1}^{s}
\left(
\displaystyle\sum_{|\alpha|\leq m}\displaystyle\sum_{i=1}^{r}
(a_{\alpha})_{ij}(p)\,D^{\alpha}u^{i}(p)
\right)\widehat{e}_{j}^{0}
\right).
\]
By definition of $a_{\alpha}(p)$, the vector in $\mathbb{K}^{s}$ whose $\,j$th component is
\[
\displaystyle\sum_{|\alpha|\leq m}\displaystyle\sum_{i=1}^{r}
(a_{\alpha})_{ij}(p)\,D^{\alpha}u^{i}(p)
\]
is precisely
\[
\displaystyle\sum_{|\alpha|\leq m}a_{\alpha}(p)\,D^{\alpha}\widetilde u(p).
\]
Thus,
\[
\pi_{2}\circ\boldsymbol{\tau}_{\mathbf{F}}\bigl(P(\mathbf{u})\bigr)(p)
=
\displaystyle\sum_{|\alpha|\leq m}a_{\alpha}(p)\,D^{\alpha}\widetilde u(p),
\qquad \forall\,p\in U.
\]

Since $(U,\phi)$ and the trivializations $\boldsymbol{\tau}_{\mathbf{E}},\boldsymbol{\tau}_{\mathbf{F}}$ were arbitrary, we have checked that, in every smooth chart and for every pair of local trivializations, there is a family of smooth maps
\[
a_{\alpha}\colon U\longrightarrow\operatorname{Hom}(\mathbb{K}^{r},\mathbb{K}^{s}),
\qquad |\alpha|\leq m,
\]
providing the expression
\[
\pi_{2}\circ\boldsymbol{\tau}_{\mathbf{F}}\bigl(P(\mathbf{u})\bigr)(p)
=
\displaystyle\sum_{|\alpha|\leq m}a_{\alpha}(p)\,D^{\alpha}
\bigl(\pi_{2}\circ\boldsymbol{\tau}_{\mathbf{E}}\circ \mathbf{u}\bigr)(p)
\]
for every $\mathbf{u}\in\Gamma(\mathbf{E})$ and every $p\in U$. By Theorem~\ref{thm:pdo expresion local explicita}, we conclude that
\[
P\in\mathbf{PDO}^{(m)}(\mathbf{E},\mathbf{F}).
\]
\end{proof}

\begin{proposition}[Conjugation of differential operators]
\label{prop:conjugacion-operadores-diferenciales-haces}
Let $M$ be a smooth manifold with or without boundary. Let $\mathbf{E},\mathbf{F}\to M$ be complex vector bundles and let
\[
P\colon\Gamma(\mathbf{E})\longrightarrow\Gamma(\mathbf{F})
\]
be a complex-linear operator. Define
\[
\overline P\colon\Gamma(\overline{\mathbf{E}})\longrightarrow
\Gamma(\overline{\mathbf{F}}),
\qquad
\overline P(\overline{\mathbf{u}}):=\overline{P\mathbf{u}}.
\]
Then:
\begin{enumerate}[label=(\alph*)]
\item $\overline P$ is well defined and complex-linear;
\item if $P\in\mathbf{PDO}^{(m)}(\mathbf{E},\mathbf{F})$, then $\overline P\in\mathbf{PDO}^{(m)}(\overline{\mathbf{E}},\overline{\mathbf{F}})$;
\item if, in a chart and with respect to local frames $(\mathbf{e}_1,\dots,\mathbf{e}_r)$ and $(\mathbf{f}_1,\dots,\mathbf{f}_s)$, the expression for $P$ is
\[
(P\mathbf{u})^{\sim}=\sum_{|\alpha|\leq m}A_\alpha D^\alpha\widetilde u,
\qquad
A_\alpha\colon U\longrightarrow
M_{s\times r}(\mathbb C),
\]
then, with respect to the conjugate frames $(\overline{\mathbf{e}}_1,\dots,\overline{\mathbf{e}}_r)$ and $(\overline{\mathbf{f}}_1,\dots,\overline{\mathbf{f}}_s)$,
\[
(\overline P\mathbf{z})^{\sim}
=\sum_{|\alpha|\leq m}\overline{A_\alpha}\,D^\alpha\widetilde z;
\]
\item conjugation commutes with restrictions and satisfies
\[
\overline{Q\circ P}=\overline Q\circ\overline P,
\qquad
\overline{\operatorname{id}_{\mathbf{E}}}=\operatorname{id}_{\overline{\mathbf{E}}}.
\]
In particular, if $P$ is invertible, then $\overline{P^{-1}}=(\overline P)^{-1}$.
\end{enumerate}
\end{proposition}

\begin{proof}
Conjugation of sections is an antilinear bijection, so every section of $\overline{\mathbf{E}}$ is uniquely expressible as $\overline{\mathbf{u}}$. If $\lambda\in\mathbb C$, the scalar structure of $\overline{\mathbf{E}}$ gives $\lambda\overline{\mathbf{u}}=\overline{\overline\lambda \mathbf{u}}$; therefore,
\[
\overline P(\lambda\overline{\mathbf{u}})
=\overline{P(\overline\lambda \mathbf{u})}
=\overline{\overline\lambda P\mathbf{u}}
=\lambda\overline{P\mathbf{u}}.
\]
Additivity is checked in the same way.

For the local formula, write $\displaystyle \mathbf{z}=\displaystyle\sum_{j=1}^{r} z^j\overline{\mathbf{e}}_j$. Then $\mathbf{z}=\overline{\mathbf{u}}$ for $\displaystyle \mathbf{u}=\displaystyle\sum_{j=1}^{r}\overline{z^j}\mathbf{e}_j$. Conjugating the local expression for $P\mathbf{u}$ gives
\[
(\overline P\mathbf{z})^{\sim}
=\overline{\sum_{|\alpha|\leq m}
A_\alpha D^\alpha(\overline{\widetilde z})}
=\sum_{|\alpha|\leq m}
\overline{A_\alpha}D^\alpha\widetilde z.
\]
The local characterization of differential operators then shows that $\overline P$ has order at most $m$. Compatibility with restrictions is immediate from the same formula. Finally, for every $\mathbf{u}\in\Gamma(\mathbf{E})$,
\[
\overline{Q\circ P}(\overline{\mathbf{u}})
=\overline{Q(P\mathbf{u})}
=(\overline Q\circ\overline P)(\overline{\mathbf{u}}),
\]
and the assertions about identities and inverses follow from this equality.
\end{proof}

\section{Formal Hermitian adjoints}

Throughout the complex case, we use Hermitian metrics linear in the first variable and conjugate-linear in the second. All integrals are taken with respect to the Riemann--Lebesgue measure $\lambda_{\mathbf g}$ constructed in Definition~\ref{medibles ajenos dos a dos, definicion de medida Riemann--Lebesgue}; in particular, no orientation of $M$ is required.

\begin{definition}[Formal Hermitian adjoint]
\label{def: adjunto formal de un pdo}\glsadd{adjunto-formal}
\label{def:adjunto-hermitiano-formal}
Let $(M,\mathbf{g})$ be a Riemannian manifold without boundary and let $\mathbf{E},\mathbf{F}\to M$ be real vector bundles with bundle metrics if $\mathbb K=\mathbb R$, or Hermitian bundles if $\mathbb K=\mathbb C$. For $P\in\mathbf{PDO}(\mathbf{E},\mathbf{F})$, an operator
\[
P_h^*\in\mathbf{PDO}(\mathbf{F},\mathbf{E})
\]
is called a \emph{formal Hermitian adjoint} of $P$ if
\[
\int_M \mathbf{h}_{\mathbf{F}}(P\mathbf{u},\mathbf{v})\,d\lambda_{\mathbf{g}}
=\int_M \mathbf{h}_{\mathbf{E}}(\mathbf{u},P_h^*\mathbf{v})\,d\lambda_{\mathbf{g}}
\]
for every $\mathbf{u}\in\Gamma_c(\mathbf{E})$ and every $\mathbf{v}\in\Gamma_c(\mathbf{F})$. In the real case, this is the usual real notion.
\end{definition}

\begin{lemma}[Uniqueness]
\label{lema: unicidad adjunto formal}
Let $(M,\mathbf{g})$ be a Riemannian manifold without boundary and let $\mathbf{E},\mathbf{F}\to M$ be real vector bundles with bundle metrics or Hermitian bundles. A differential operator $P\in\mathbf{PDO}(\mathbf{E},\mathbf{F})$ admits at most one formal Hermitian adjoint.
\end{lemma}

\begin{proof}
Let $S_1,S_2\in\mathbf{PDO}(\mathbf{F},\mathbf{E})$ be two adjoints and fix $\mathbf{v}\in\Gamma_c(\mathbf{F})$. Locality implies that $\mathbf{w}:=(S_1-S_2)\mathbf{v}$ belongs to $\Gamma_c(\mathbf{E})$. Taking $\mathbf{u}=\mathbf{w}$ in the difference of the two adjunction identities gives
\[
0=\int_M \mathbf{h}_{\mathbf{E}}(\mathbf{w},\mathbf{w})\,d\lambda_{\mathbf{g}},
\]
and, by positivity, $\mathbf{w}=0$. Thus, $S_1\mathbf{v}=S_2\mathbf{v}$ for every compactly supported section. Now given $\mathbf{v}\in\Gamma(\mathbf{F})$ and $x\in M$, choose $\chi\in C_c^\infty(M)$ equal to one on a neighborhood of $x$. By locality, $(S_1-S_2)\mathbf{v}=(S_1-S_2)(\chi \mathbf{v})$ near $x$, and the last expression is zero. Since $x$ was arbitrary, $S_1=S_2$.
\end{proof}

\begin{theorem}[Existence and local formula for the Hermitian adjoint]
\label{teo: adjunto pdo trivial}
\label{cor: adjunto global en coordenadas locales (existencia y expresion)}
\label{teo:existencia-adjunto-hermitiano-formal}
Let $(M,\mathbf{g})$ be a Riemannian manifold without boundary and let $\mathbf{E},\mathbf{F}\to M$ be real vector bundles with bundle metrics or Hermitian bundles. Let $P\in\mathbf{PDO}^{(m)}(\mathbf{E},\mathbf{F})$. Then there exists a unique $P_h^*\in\mathbf{PDO}^{(m)}(\mathbf{F},\mathbf{E})$.

More precisely, let $\phi\colon U\to\Omega\subseteq\mathbb R^n$ be a chart, with coordinate functions $(x^1,\dots,x^n)$, and choose orthonormal local frames (unitary in the complex case) $\mathbf e_1,\dots,\mathbf e_r$ and $\mathbf f_1,\dots,\mathbf f_s$ for $\mathbf E$ and $\mathbf F$, respectively. For a scalar function $a\in C^\infty(U,\mathbb K)$ write
\[
D_\phi^\alpha a
:=
\bigl[D^\alpha(a\circ\phi^{-1})\bigr]\circ\phi.
\]
If the functions $a_{\alpha i}^{\,j}\in C^\infty(U,\mathbb K)$ are determined by
\begin{equation}
\label{eq:expresion-local-componentes-adjunto-hermitiano}
P\left(\sum_{i=1}^r u^i\mathbf e_i\right)
=
\sum_{j=1}^s
\left(
\sum_{i=1}^r\sum_{|\alpha|\leq m}
a_{\alpha i}^{\,j}D_\phi^\alpha u^i
\right)\mathbf f_j,
\end{equation}
and $g_{ab}:=\mathbf g(\boldsymbol{\partial}_a,
\boldsymbol{\partial}_b)$ are the coordinate coefficients of the metric, then the expression for $P_h^*$ with respect to the same frames is
\begin{equation}
\label{eq:formula-local-adjunto-hermitiano}
P_h^*\left(\sum_{j=1}^s v^j\mathbf f_j\right)
=
\frac{1}{\sqrt{\det(\mathbf{g})}}
\sum_{i=1}^r
\left[
\sum_{j=1}^s\sum_{|\alpha|\leq m}(-1)^{|\alpha|}
D_\phi^\alpha\!\left(
\sqrt{\det(\mathbf{g})}\,
\overline{a_{\alpha i}^{\,j}}\,v^j
\right)
\right]\mathbf e_i.
\end{equation}
In the real case, the bar does not alter the coefficient. If the coefficients are grouped into the matrix-valued map
\[
A_\alpha=(a_{\alpha i}^{\,j})\colon
U\longrightarrow M_{s\times r}(\mathbb K),
\]
its pointwise conjugate transpose is denoted by $A_\alpha^*=\overline{A_\alpha}^{\,T}$.
\end{theorem}

\begin{proof}
Fix a chart $\phi\colon U\to\Omega$ and orthonormal local frames $\mathbf e_1,\dots,\mathbf e_r$ and $\mathbf f_1,\dots,\mathbf f_s$ for $\mathbf E$ and $\mathbf F$, respectively. The expression on the right-hand side of \eqref{eq:formula-local-adjunto-hermitiano} defines an operator
\[
Q_U\colon\Gamma(\mathbf F|_U)\longrightarrow\Gamma(\mathbf E|_U).
\]
Indeed, $\det(\mathbf{g})$ is smooth and strictly positive. Expanding each term using the Leibniz rule for multi-indices, $Q_U\mathbf v$ depends on the coordinate derivatives of the functions $v^j$ of order at most $m$. The local characterization in Corollary~\ref{cor:pdo-expresion-con-marcos} therefore implies that $Q_U\in\mathbf{PDO}^{(m)}(\mathbf F|_U,\mathbf E|_U)$.

We prove the local adjunction identity. Let $\mathbf u\in\Gamma_c(\mathbf E|_U)$ and let $\mathbf v\in\Gamma(\mathbf F|_U)$, and write
\[
\mathbf u=\sum_{i=1}^r u^i\mathbf e_i,
\qquad
\mathbf v=\sum_{j=1}^s v^j\mathbf f_j.
\]
The image $\phi(\operatorname{supp}\mathbf u)$ is a compact set contained in $\Omega$. Since the frames are orthonormal and the Hermitian metric is linear in the first variable, the local formula for the Riemann--Lebesgue measure and \eqref{eq:expresion-local-componentes-adjunto-hermitiano} give
\begin{align*}
&\int_U\mathbf h_{\mathbf F}(P\mathbf u,\mathbf v)
\,d\lambda_{\mathbf g}\\
&\quad=
\sum_{i=1}^r\sum_{j=1}^s\sum_{|\alpha|\leq m}
\int_\Omega
\bigl(a_{\alpha i}^{\,j}\circ\phi^{-1}\bigr)
D^\alpha\bigl(u^i\circ\phi^{-1}\bigr)
\overline{\bigl(v^j\circ\phi^{-1}\bigr)}
\bigl(\sqrt{\det(\mathbf{g})}\circ\phi^{-1}\bigr)
\,d\lambda_n.
\end{align*}

Fix $i$, $j$, and $\alpha=(\alpha_1,\dots,\alpha_n)$. We can integrate $\alpha_a$ times with respect to $x^a$, for each $a$, because $u^i\circ\phi^{-1}$ has compact support in $\Omega$. No boundary term appears, and each integration introduces a minus sign. Therefore,
\begin{align*}
&\int_\Omega
\bigl(a_{\alpha i}^{\,j}\circ\phi^{-1}\bigr)
D^\alpha\bigl(u^i\circ\phi^{-1}\bigr)
\overline{\bigl(v^j\circ\phi^{-1}\bigr)}
\bigl(\sqrt{\det(\mathbf{g})}\circ\phi^{-1}\bigr)
\,d\lambda_n\\
&\quad=(-1)^{|\alpha|}
\int_\Omega
\bigl(u^i\circ\phi^{-1}\bigr)
D^\alpha\!\left[
\left(
a_{\alpha i}^{\,j}\,\overline{v^j}\,
\sqrt{\det(\mathbf{g})}
\right)\circ\phi^{-1}
\right]d\lambda_n.
\end{align*}
All coordinate derivatives are real and $\sqrt{\det(\mathbf{g})}$ is real. Consequently,
\[
D^\alpha\!\left[
\left(
a_{\alpha i}^{\,j}\,\overline{v^j}\,
\sqrt{\det(\mathbf{g})}
\right)\circ\phi^{-1}
\right]
=
\overline{
D^\alpha\!\left[
\left(
\overline{a_{\alpha i}^{\,j}}\,v^j\,
\sqrt{\det(\mathbf{g})}
\right)\circ\phi^{-1}
\right]}.
\]
This equality includes the real case, in which conjugation is the identity. Substituting it, summing over $i$, $j$, and $\alpha$, and inserting the factor $\bigl(\sqrt{\det(\mathbf{g})}\bigr)^{-1}
\sqrt{\det(\mathbf{g})}$ gives
\begin{align*}
&\int_U\mathbf h_{\mathbf F}(P\mathbf u,\mathbf v)
\,d\lambda_{\mathbf g}\\
&\quad=
\int_\Omega\sum_{i=1}^r
\bigl(u^i\circ\phi^{-1}\bigr)
\overline{\left[
\frac{1}{\sqrt{\det(\mathbf{g})}}
\sum_{j=1}^s\sum_{|\alpha|\leq m}(-1)^{|\alpha|}
D_\phi^\alpha\!\left(
\sqrt{\det(\mathbf{g})}\,
\overline{a_{\alpha i}^{\,j}}\,v^j
\right)
\right]\circ\phi^{-1}}
\bigl(\sqrt{\det(\mathbf{g})}\circ\phi^{-1}\bigr)
\,d\lambda_n\\
&\quad=
\int_U\mathbf h_{\mathbf E}(\mathbf u,Q_U\mathbf v)
\,d\lambda_{\mathbf g}.
\end{align*}
We have proved a slightly stronger identity than required: only $\mathbf u$ needs to have compact support in $U$; the section $\mathbf v$ may be any smooth section over $U$.

We now show that the local operators are compatible. Let $U$ and $V$ be two coordinate domains with orthonormal frames as above, and let $W=U\cap V$. For $\mathbf u\in\Gamma_c(\mathbf E|_W)$ and $\mathbf v\in\Gamma_c(\mathbf F|_W)$, the two formulas already proved give
\[
\int_W\mathbf h_{\mathbf E}
\bigl(\mathbf u,(Q_U|_W)\mathbf v\bigr)d\lambda_{\mathbf g}
=\int_W\mathbf h_{\mathbf F}(P\mathbf u,\mathbf v)d\lambda_{\mathbf g}
=\int_W\mathbf h_{\mathbf E}
\bigl(\mathbf u,(Q_V|_W)\mathbf v\bigr)d\lambda_{\mathbf g}.
\]
Thus, $Q_U|_W$ and $Q_V|_W$ are formal adjoints of the same operator $P|_W$. Lemma~\ref{lema: unicidad adjunto formal}, applied to the Riemannian manifold $W$, implies that $Q_U|_W=Q_V|_W$. In particular, the equality is not assumed from the coordinate expressions; it follows from the integral identity both satisfy.

We can cover $M$ by a locally finite family $\{U_\ell\}_{\ell\in L}$ of coordinate domains over which $\mathbf E$ and $\mathbf F$ admit orthonormal frames. Compatibility on intersections allows us to define, for $\mathbf v\in\Gamma(\mathbf F)$, a section $P_h^*\mathbf v$ by
\[
(P_h^*\mathbf v)|_{U_\ell}:=Q_{U_\ell}(\mathbf v|_{U_\ell}).
\]
The sections on the right-hand side agree on overlaps, so they define a unique global smooth section. Since the order of a differential operator can be checked locally, the formulas for the $Q_{U_\ell}$ show that
\[
P_h^*\in\mathbf{PDO}^{(m)}(\mathbf F,\mathbf E).
\]

It remains to check the global identity. Let $\{\chi_\ell\}_{\ell\in L}$ be a smooth partition of unity subordinate to $\{U_\ell\}_{\ell\in L}$, and let $\mathbf u\in\Gamma_c(\mathbf E)$ and $\mathbf v\in\Gamma_c(\mathbf F)$. The family $\{\operatorname{supp}(\chi_\ell)\}_{\ell\in L}$ is locally finite. Thus each point of the compact set $\operatorname{supp}(\mathbf u)$ has a neighborhood meeting only finitely many of those supports; a finite subcover by such neighborhoods shows that the set
\[
L_{\mathbf u}:=
\{\ell\in L\mid
\operatorname{supp}(\mathbf u)\cap
\operatorname{supp}(\chi_\ell)\neq\varnothing\}
\]
is finite. Consequently,
\[
\mathbf u=\sum_{\ell\in L_{\mathbf u}}\chi_\ell\mathbf u.
\]
Moreover, $\chi_\ell\mathbf u$ has compact support contained in $U_\ell$. Locality of $P$ implies
\[
\operatorname{supp}\bigl(P(\chi_\ell\mathbf u)\bigr)
\subseteq\operatorname{supp}(\chi_\ell\mathbf u)
\subseteq U_\ell.
\]
We can apply on $U_\ell$ the local identity in which the first section has compact support and the second is $\mathbf v|_{U_\ell}$. Thus,
\begin{align*}
\int_M\mathbf h_{\mathbf F}(P\mathbf u,\mathbf v)d\lambda_{\mathbf g}
&=\sum_{\ell\in L_{\mathbf u}}
\int_{U_\ell}\mathbf h_{\mathbf F}
\bigl(P(\chi_\ell\mathbf u),\mathbf v\bigr)d\lambda_{\mathbf g}\\
&=\sum_{\ell\in L_{\mathbf u}}
\int_{U_\ell}\mathbf h_{\mathbf E}
\bigl(\chi_\ell\mathbf u,P_h^*\mathbf v\bigr)d\lambda_{\mathbf g}\\
&=\int_M\mathbf h_{\mathbf E}(\mathbf u,P_h^*\mathbf v)
d\lambda_{\mathbf g}.
\end{align*}
Therefore, $P_h^*$ is a formal Hermitian adjoint of $P$. Its uniqueness is that of Lemma~\ref{lema: unicidad adjunto formal}.
\end{proof}

\begin{remark}[Formal adjoint in the presence of a boundary]
\label{obs:adjunto-formal-frontera-dominio}
If $M$ has smooth boundary and the coefficients of $P$, the bundle metrics, and the base metric are smooth up to it, the same local formula defines a differential operator $P_h^*$ smooth up to the boundary. Indeed, in a boundary chart the determinant of the metric remains positive, and all derivatives of the coefficients in \eqref{eq:formula-local-adjunto-hermitiano} are smooth up to the flat part. On overlaps, the expressions agree first in the interior by the uniqueness already proved, and then on the boundary by continuity. The identity characterizing this formal adjoint is required for test sections with compact support contained in $\operatorname{Int}(M)$.

When the sections reach the boundary, integration by parts may produce additional terms, such as those we will calculate for covariant derivatives. The formal differential operator is determined before boundary conditions are chosen. Its realization as an unbounded operator on $L^2(M,\mathbf E)$, and the domain of its Hilbert space adjoint, require specifying a domain of sections; they are not identified by the differential formula alone.
\end{remark}

\begin{proposition}[Rules for the formal Hermitian adjoint]
\label{adjunto formal composición}
\label{prop:reglas-adjunto-hermitiano-formal}
Let $(M,\mathbf{g})$ be a Riemannian manifold without boundary, let $\mathbf{E},\mathbf{F},\mathbf{G}\to M$ be real vector bundles with bundle metrics or Hermitian bundles, and let $P\in\mathbf{PDO}(\mathbf{E},\mathbf{F})$ and $Q\in\mathbf{PDO}(\mathbf{F},\mathbf{G})$. Then:
\begin{enumerate}[label=(\alph*)]
\item $(Q\circ P)_h^*=P_h^*\circ Q_h^*$;
\item for $R\in\mathbf{PDO}(\mathbf{E},\mathbf{F})$ and $\lambda\in\mathbb K$,
\[
(P+\lambda R)_h^*=P_h^*+\overline\lambda\,R_h^*,
\]
where $\overline\lambda=\lambda$ in the real case;
\item if $M_f$ denotes multiplication by $f\in C^\infty(M,\mathbb K)$, then $(M_f)_h^*=M_{\overline f}$;
\item for every open set $U\subseteq M$,
\[
(P|_U)_{h|_U}^*=(P_h^*)|_U;
\]
\item if $\mathbf{E},\mathbf{F}$ are complex and $\overline{\mathbf{E}},\overline{\mathbf{F}}$ are equipped with the conjugate metrics, then
\[
(\overline P)_{\overline h}^*=\overline{P_h^*}.
\]
\end{enumerate}
\end{proposition}

\begin{proof}
The first identity follows by successively transferring $Q$ and $P$ in the integral; supports remain compact by locality. For the second, linearity in the first variable and conjugate-linearity in the second give
\[
\mathbf{h}_{\mathbf{E}}(\mathbf{u},\overline\lambda R_h^*\mathbf{v})
=\lambda \mathbf{h}_{\mathbf{E}}(\mathbf{u},R_h^*\mathbf{v}).
\]
The third is proved in the same way. The fourth follows by applying the integral identity to compactly supported sections in $U$ and using uniqueness.

For the last assertion, write the conjugate metrics as $\mathbf{h}_{\overline{\mathbf{E}}}(\overline{\mathbf{u}},\overline{\mathbf{w}})
=\overline{\mathbf{h}_{\mathbf{E}}(\mathbf{u},\mathbf{w})}$ and $\mathbf{h}_{\overline{\mathbf{F}}}(\overline{\mathbf{v}},\overline{\mathbf{z}})
=\overline{\mathbf{h}_{\mathbf{F}}(\mathbf{v},\mathbf{z})}$. Then
\begin{align*}
\int_M \mathbf{h}_{\overline{\mathbf{F}}}(\overline P\,\overline{\mathbf{u}},\overline{\mathbf{v}})
\,d\lambda_{\mathbf{g}}
&=\overline{\int_M \mathbf{h}_{\mathbf{F}}(P\mathbf{u},\mathbf{v})\,d\lambda_{\mathbf{g}}}\\
&=\overline{\int_M \mathbf{h}_{\mathbf{E}}(\mathbf{u},P_h^*\mathbf{v})\,d\lambda_{\mathbf{g}}}\\
&=\int_M \mathbf{h}_{\overline{\mathbf{E}}}
(\overline{\mathbf{u}},\overline{P_h^*}\,\overline{\mathbf{v}})\,d\lambda_{\mathbf{g}},
\end{align*}
and uniqueness completes the proof.
\end{proof}

\begin{exercise}
\label{linealidad del adjunto formal}
Check the rule $(P+\lambda R)_h^*=P_h^*+\overline\lambda R_h^*$ directly and explain why the bar is indispensable when $\mathbb K=\mathbb C$.
\end{exercise}

\begin{definition}[Formally self-adjoint operator]
\label{def:operadores-diferenciales-en-haces-formalmente-autoadjunto}
Let $(M,\mathbf{g})$ be a Riemannian manifold without boundary and let $\mathbf{E}\to M$ be a real vector bundle with a bundle metric or a Hermitian bundle. An operator $P\in\mathbf{PDO}(\mathbf{E},\mathbf{E})$ is \emph{formally self-adjoint} if $P=P_h^*$.
\end{definition}
We conclude this section with examples illustrating the concept of a formal Hermitian adjoint. We begin with the directional covariant derivative.

\begin{proposition}[Formal adjoint of $\nabla_{\mathbf{X}}$]\label{prop:adjunto-formal-nablaX}\index{formal adjoint}
Let $(M,\mathbf{g})$ be a Riemannian manifold (with or without boundary $\partial M$) and let $\nabla$ be the Levi--Civita connection. Let $\mathbf{X}\in \mathfrak{X}(M)$ be a fixed smooth vector field. Then the operator
\[
\nabla_{\mathbf{X}}\colon \mathfrak{X}(M)\longrightarrow \mathfrak{X}(M),
\qquad
\mathbf{Y}\longmapsto \nabla_{\mathbf{X}} \mathbf{Y},
\]
is a partial differential operator of order $1$, that is,
\[
\nabla_{\mathbf{X}}\in \mathbf{PDO}^{(1)}(TM,TM).
\]

Moreover, for any vector fields $\mathbf{Y},\mathbf{Z}\in \mathfrak{X}_c(M)$ we have
\begin{equation}\label{eq:adjunto-nablaX-TM}
\int_M \mathbf{g}(\nabla_{\mathbf{X}} \mathbf{Y}, \mathbf{Z})\,d\lambda_{\mathbf{g}}
=
-\int_M \mathbf{g}\bigl(\mathbf{Y},\nabla_{\mathbf{X}} \mathbf{Z}\bigr)\,d\lambda_{\mathbf{g}}
-\int_M (\operatorname{div}\mathbf{X})\,\mathbf{g}(\mathbf{Y},\mathbf{Z})\,d\lambda_{\mathbf{g}}
+
\int_{\partial M} \mathbf{g}(\mathbf{Y},\mathbf{Z})\,\mathbf{g}(\mathbf{X},\boldsymbol{\nu})\,d\lambda_{\widetilde{\mathbf{g}}},
\end{equation}
where, if $\partial M\neq\varnothing$, $\boldsymbol{\nu}$ is the outward unit normal to $\partial M$.
\end{proposition}

\begin{proof}
First we check that $\nabla_{\mathbf{X}}$ is a partial differential operator of order $1$ in the commutator sense.

Let $f\in C^\infty(M)$ and $\mathbf{Y}\in \mathfrak{X}(M)$. Using the Leibniz rule for the connection gives
\[
\nabla_{\mathbf{X}}(f \mathbf{Y})
=
\mathbf{X}(f)\,\mathbf{Y}+f\nabla_{\mathbf{X}} \mathbf{Y}.
\]
Therefore,
\[
[\nabla_{\mathbf{X}},f]\mathbf{Y}
=
\nabla_{\mathbf{X}}(f \mathbf{Y})-f\nabla_{\mathbf{X}} \mathbf{Y}
=
\mathbf{X}(f)\mathbf{Y}.
\]
This depends only on $\mathbf{Y}$ and not on its derivatives, so the operator $[\nabla_{\mathbf{X}},f]$ is $C^\infty(M)$-linear. Consequently, $\nabla_{\mathbf{X}}\in \mathbf{PDO}^{(1)}(TM,TM)$.

We now prove the integration by parts identity. Let $\mathbf{Y},\mathbf{Z}\in \mathfrak{X}_c(M)$ and define
\[
f:=\mathbf{g}(\mathbf{Y},\mathbf{Z}),\qquad \mathbf{V}:=f\mathbf{X}.
\]
By the Leibniz rule for divergence (Proposition~\ref{leibniz divergencia}),
\[
\operatorname{div}(f\mathbf{X})=\mathbf{X}(f)+f\,\operatorname{div}\mathbf{X}.
\]
Since the Levi--Civita connection is compatible with the metric,
\[
\mathbf{X}(f)=\mathbf{X}\bigl(\mathbf{g}(\mathbf{Y},\mathbf{Z})\bigr)=\mathbf{g}(\nabla_{\mathbf{X}} \mathbf{Y},\mathbf{Z})+\mathbf{g}(\mathbf{Y},\nabla_{\mathbf{X}} \mathbf{Z}).
\]
Consequently,
\[
\operatorname{div}(f\mathbf{X})
=
\mathbf{g}(\nabla_{\mathbf{X}} \mathbf{Y},\mathbf{Z})+\mathbf{g}(\mathbf{Y},\nabla_{\mathbf{X}} \mathbf{Z})+(\operatorname{div}\mathbf{X})\mathbf{g}(\mathbf{Y},\mathbf{Z}).
\]
Integrate over $M$ and apply the divergence theorem:
\[
\int_M \operatorname{div}(f\mathbf{X})\,d\lambda_{\mathbf{g}}
=
\int_{\partial M} \mathbf{g}(f\mathbf{X},\boldsymbol{\nu})\,d\lambda_{\widetilde{\mathbf{g}}}
=
\int_{\partial M} \mathbf{g}(\mathbf{Y},\mathbf{Z})\,\mathbf{g}(\mathbf{X},\boldsymbol{\nu})\,d\lambda_{\widetilde{\mathbf{g}}}.
\]
Rearranging the terms gives \eqref{eq:adjunto-nablaX-TM}.
\end{proof}

\begin{remark}[Formal adjoint of $\nabla_{\mathbf{X}}$]\label{obs:operadores-diferenciales-en-haces-adjunto-formal-de}
Let $(M,\mathbf{g})$ be a Riemannian manifold without boundary and let $\mathbf{X}\in\mathfrak X(M)$. For $\mathbf{Y},\mathbf{Z}\in\mathfrak X_c(M)$, identity \eqref{eq:adjunto-nablaX-TM} takes the form
\[
\int_M \mathbf{g}(\nabla_{\mathbf{X}} \mathbf{Y}, \mathbf{Z})\,d\lambda_{\mathbf{g}}
=
\int_M \mathbf{g}\bigl(\mathbf{Y},-\nabla_{\mathbf{X}} \mathbf{Z}-(\operatorname{div}\mathbf{X})\mathbf{Z}\bigr)\,d\lambda_{\mathbf{g}}.
\]
In this sense, the formal adjoint of the partial differential operator
\[
\nabla_{\mathbf{X}}\colon \mathfrak{X}(M)\longrightarrow \mathfrak{X}(M)
\]
is given by
\[
\nabla_{\mathbf{X}}^{*}\mathbf{Z}=-\nabla_{\mathbf{X}} \mathbf{Z}-(\operatorname{div}\mathbf{X})\,\mathbf{Z}.
\]

This operator should be distinguished from the global formal adjoint
\[
\nabla^{*}\colon \Gamma(T^*M\otimes TM)\longrightarrow \mathfrak{X}(M),
\]
corresponding to the covariant derivative viewed as a first-order differential operator between vector bundles. Later we will derive an integration by parts formula allowing us to compute $\nabla^{*}$ explicitly.
\end{remark}

Having introduced the formal adjoint of the covariant derivative, we can define a second-order operator fundamental to Riemannian geometry.

\begin{definition}[Bochner Laplacian]\label{def:operadores-diferenciales-en-haces-laplaciano-de-bochner}\index{Bochner Laplacian}\index{differential operator!Bochner Laplacian}\glsadd{laplaciano-bochner}
Let $(M,\mathbf{g})$ be a Riemannian manifold without boundary and let $\pi_{\mathbf{E}}\colon \mathbf{E}\longrightarrow M$ be a smooth vector bundle equipped with a bundle metric $\mathbf{h}_{\mathbf{E}}$ and a compatible connection $\nabla^{\mathbf{E}}$. We define the Bochner Laplacian of a section $\mathbf{u}\in \Gamma(\mathbf{E})$ by
\[
\Delta_B \mathbf{u}
:=
\nabla^{\mathbf{E}*}\nabla^{\mathbf{E}} \mathbf{u},
\]
where $\nabla^{\mathbf{E}*}$ denotes the formal adjoint of
\[
\nabla^{\mathbf{E}}\colon \Gamma(\mathbf{E})\longrightarrow \Gamma(T^*M\otimes \mathbf{E}).
\]
\end{definition}

\begin{remark}\label{obs:operadores-diferenciales-en-haces-laplaciano-bochner-formalmente-autoadjunto-efecto}
The Bochner Laplacian is formally self-adjoint, as shown by the identity
\[
(\nabla^{\mathbf{E}*}\nabla^{\mathbf{E}})^{*}
=
(\nabla^{\mathbf{E}})^{*}(\nabla^{\mathbf{E}*})^{*}
=
\nabla^{\mathbf{E}*}\nabla^{\mathbf{E}},
\]
since $(\nabla^{\mathbf{E}*})^{*}=\nabla^{\mathbf{E}}$.

\end{remark}
Besides the Bochner Laplacian, there are other important Laplacians, such as the Hodge--de Rham Laplacian, defined on spaces of differential forms:
\begin{definition}[Space of differential forms]\label{def:operadores-diferenciales-en-haces-espacio-de-formas-diferenciales}\index{space of differential forms}
Let $M$ be a smooth manifold with or without boundary of dimension $n$. Recall that the space of differential $k$-forms is defined by
\[
\Omega^{k}(M):=\Gamma\bigl(\Lambda^{k}T^{*}M\bigr),
\qquad k\in\{0,\dots,n\}.
\]
We also define the graded space of all forms as the finite direct sum of $C^\infty(M)$-modules
\[
\Omega^{\bullet}(M):=\bigoplus_{k=0}^{n}\Omega^{k}(M).
\]
\end{definition}
Now let $(M,\mathbf{g})$ be a Riemannian manifold without boundary. We use the determinant inner product of
Example~\ref{ej:modulo-clifford-formas} on each exterior bundle $\Lambda^{k}T^{*}M$, denoted by $\langle\cdot,\cdot\rangle_{\mathbf{g}}$, and the Riemannian measure $d\lambda_{\mathbf{g}}$ allows us to define the global inner product
\[
\langle \boldsymbol{\alpha},\boldsymbol{\beta}\rangle_{L^{2}(M,\Lambda^kT^*M)}
:=
\int_{M}\langle \boldsymbol{\alpha},\boldsymbol{\beta}\rangle_{\mathbf{g}}\,d\lambda_{\mathbf{g}},
\qquad
\boldsymbol{\alpha},\boldsymbol{\beta}\in \Omega^{k}_c(M).
\]

Recall that the exterior derivative
\[
d\colon \Omega^{k}(M)\longrightarrow \Omega^{k+1}(M)
\]
is a partial differential operator of order $1$.

If $f\in C^{\infty}(M)$ and $\boldsymbol{\alpha}\in \Omega^{k}(M)$, the Leibniz rule gives
\[
d(f\boldsymbol{\alpha})=df\wedge \boldsymbol{\alpha}+fd\boldsymbol{\alpha},
\]
and therefore
\[
[d,f]\boldsymbol{\alpha}=d(f\boldsymbol{\alpha})-fd\boldsymbol{\alpha}
=
df\wedge \boldsymbol{\alpha}.
\]
This depends only on $\boldsymbol{\alpha}$ and not on its derivatives, so the operator $[d,f]$ is $C^{\infty}(M)$-linear. Consequently,
\[
d\in \mathbf{PDO}^{(1)}\bigl(\Lambda^{k}T^{*}M,\Lambda^{k+1}T^{*}M\bigr).
\]

\begin{definition}[Codifferential]\label{def:codiferencial}\index{codifferential}
Let $(M,\mathbf{g})$ be a Riemannian manifold without boundary of dimension $n$. The formal adjoint of the exterior derivative
\[
d\colon \Omega^{k}(M)\longrightarrow \Omega^{k+1}(M)
\]
is denoted by
\[
\delta\colon \Omega^{k+1}(M)\longrightarrow \Omega^{k}(M)
\]
and is called the codifferential.

If $M$ is oriented, the Hodge star operator of
Definition~\ref{def:operadores-diferenciales-en-haces-operador-de-hodge} gives the following expression on
$\Omega^{k+1}(M)$, proved in Proposition~\ref{prop:operadores-diferenciales-en-haces-orientable}:
\[
\delta
=
(-1)^{k+1}\,*^{-1}d*.
\]
\end{definition}

This formal adjoint is determined by the rule
\[
\int_M \langle d\boldsymbol{\alpha},\boldsymbol{\beta}\rangle_{\mathbf{g}}\,d\lambda_{\mathbf{g}}
=
\int_M \langle \boldsymbol{\alpha},\delta\boldsymbol{\beta}\rangle_{\mathbf{g}}\,d\lambda_{\mathbf{g}},
\]
for $\boldsymbol{\alpha}\in\Omega_{c}^{k}(M)$ and $\boldsymbol{\beta}\in\Omega_{c}^{k+1}(M)$.
\begin{definition}[Hodge--de Rham Laplacian]\label{def:operadores-diferenciales-en-haces-laplaciano-de-hodge-derham}\index{Hodge--de Rham Laplacian}\index{differential operator!Hodge--de Rham Laplacian}\glsadd{laplaciano-hodge}
Let $(M,\mathbf{g})$ be a Riemannian manifold without boundary. We define the Hodge--de Rham Laplacian to be the second-order operator
\[
\Delta
:=
d\delta+\delta d
\colon \Omega^{\bullet}(M)\longrightarrow \Omega^{\bullet}(M).
\]
\end{definition}

\begin{proposition}\label{prop:operadores-diferenciales-en-haces-laplaciano-hodge-derham-formalmente-autoadjunto}
Let $(M,\mathbf{g})$ be an oriented Riemannian manifold without boundary. The Hodge--de Rham Laplacian $\Delta$ is formally self-adjoint.
\end{proposition}

\begin{proof}
We have
\[
(d\delta)^{*}=\delta^{*}d^{*}=d\delta,
\qquad
(\delta d)^{*}=d^{*}\delta^{*}=\delta d,
\]
since $\delta=d^{*}$ and $(d^{*})^{*}=d$. Therefore,
\[
\Delta^{*}=(d\delta+\delta d)^{*}=d\delta+\delta d=\Delta.
\]
\end{proof}
The classical exterior differential operator extends naturally to the setting of vector bundles.

Let $M$ be a smooth manifold with or without boundary of dimension $n$. The space of smooth differential $k$-forms is defined by
\[
\Omega^{k}(M)
:=
\Gamma\!\left(\Lambda^{k}T^{*}M\right),
\]
that is, the space of smooth sections of the bundle $\Lambda^{k}T^{*}M$. The full space of forms is defined as the direct sum of $C^{\infty}(M)$-modules
\[
\Omega^{\bullet}(M)
:=
\bigoplus_{k=0}^{n}\Omega^{k}(M).
\]

If $(x^{1},\dots,x^{n})$ are local coordinates, we use the notation of increasing multi-indices
\[
I=(i_{1},\dots,i_{k}),
\qquad
1\le i_{1}<\dots<i_{k}\le n,
\]
and write
\[
\mathbf{d}x^{I}:=\mathbf{d}x^{i_{1}}\wedge\dots\wedge \mathbf{d}x^{i_{k}}.
\]
With this notation, a $k$-form can locally be expressed as
\[
\boldsymbol{\omega}
=
\displaystyle\sum_{1\leq i_{1}<\dots<i_{k}\leq n}
\omega_{I}\,\mathbf{d}x^{I},
\]
and the classical exterior derivative takes the form
\[
d\!\left(\displaystyle\sum_{1\leq i_{1}<\dots<i_{k}\leq n}\omega_{I}\,\mathbf{d}x^{I}\right)
=
\displaystyle\sum_{1\leq i_{1}<\dots<i_{k}\leq n} d\omega_{I}\wedge \mathbf{d}x^{I}.
\]

Now let $\pi_{\mathbf{E}}\colon \mathbf{E}\longrightarrow M$ be a smooth vector bundle of rank $r$. This allows us to generalize the concept of a differential form by considering forms with values in $\mathbf{E}$.

\begin{definition}[Bundle-valued differential forms]\label{def:operadores-diferenciales-en-haces-formas-diferenciales-con-valores-en-un-haz}\index{bundle-valued differential forms}
Let $M$ be a smooth manifold with or without boundary, let $\mathbf{E}\to M$ be a smooth vector bundle, and let $k\ge 0$. Define
\[
\Omega^{k}(M,\mathbf{E})
:=
\Gamma\!\left(\Lambda^{k}T^{*}M\otimes \mathbf{E}\right).
\]
The full space is defined as the direct sum of $C^{\infty}(M)$-modules
\[
\Omega^{\bullet}(M,\mathbf{E})
:=
\bigoplus_{k=0}^{n}\Omega^{k}(M,\mathbf{E}).
\]
\end{definition}

Observe that if $\mathbf{E}=M\times\mathbb{R}$ is the trivial bundle, the classical case is recovered:
\[
\Omega^{k}(M,\mathbf{E})\cong \Omega^{k}(M).
\]

To define an exterior derivative in this setting, a connection on the bundle $\mathbf{E}$ is needed.

\begin{definition}[Induced connection]\label{def:operadores-diferenciales-en-haces-conexion-inducida}\index{induced connection}
Let $M$ be a smooth manifold with or without boundary, let $\mathbf{E}\to M$ be a smooth vector bundle, and let $\nabla^{\mathbf{E}}$ be a connection on $\mathbf{E}$. Also fix an auxiliary connection $\nabla^M$ on $TM$, and use the same notation for its induced connections on the tensor bundles constructed from $TM$ and $T^{*}M$. The product connection $\nabla$ on $\Lambda^{k}T^{*}M\otimes \mathbf{E}$ is defined for each vector field $\mathbf{X} \in \mathfrak{X}(M)$ by the Leibniz rule
\[
\nabla_{\mathbf{X}}(\boldsymbol{\alpha}\otimes \mathbf{s})
=
\nabla^M_{\mathbf{X}}\boldsymbol{\alpha}\otimes \mathbf{s}
+
\boldsymbol{\alpha}\otimes \nabla^{\mathbf{E}}_{\mathbf{X}}\mathbf{s},
\]
for $\boldsymbol{\alpha}\in\Omega^{k}(M)$ and $\mathbf{s}\in\Gamma(\mathbf{E})$.
\end{definition}

The alternation operator will relate the product connection to the exterior derivative defined next.

\begin{definition}[Alternation operator]\label{def:operadores-diferenciales-en-haces-operador-de-alternacion}\index{alternation operator}
Let $M$ be a smooth manifold with or without boundary, let $\mathbf{E}\to M$ be a smooth vector bundle, and let
\[
\mathbf{T}\in \Gamma\left(T^{(0,k+1)}(TM)\otimes \mathbf{E}\right).
\]
Define
\[
\operatorname{Alt}(\mathbf{T})\in \Gamma\left(\Lambda^{k+1}T^{*}M\otimes \mathbf{E}\right)
\]
by
\[
\operatorname{Alt}(\mathbf{T})(\mathbf{X}_{1},\dots,\mathbf{X}_{k+1})
:=
\frac{1}{(k+1)!}
\displaystyle\sum_{\sigma\in S_{k+1}}
\operatorname{sgn}(\sigma)\,
\mathbf{T}(\mathbf{X}_{\sigma(1)},\dots,\mathbf{X}_{\sigma(k+1)}),
\]
for vector fields $\mathbf{X}_{1},\dots,\mathbf{X}_{k+1}\in \mathfrak{X}(M)$.
\end{definition}

\begin{definition}[Bundle-valued exterior derivative]\label{def:operadores-diferenciales-en-haces-derivada-exterior-con-valores-en-un-haz}\index{bundle-valued exterior derivative}
Let $M$ be a smooth manifold with or without boundary, let $\mathbf{E}\to M$ be a smooth vector bundle over $\mathbb K$, and let $\nabla^{\mathbf{E}}$ be a connection on $\mathbf{E}$. For each $k\geq0$ define
\[
d^{\mathbf{E}}\colon \Omega^{k}(M,\mathbf{E})\longrightarrow \Omega^{k+1}(M,\mathbf{E})
\]
to be the $\mathbb K$-linear operator determined locally by the rule
\[
d^{\mathbf{E}}(\boldsymbol{\alpha}\otimes\mathbf{s})
=
d\boldsymbol{\alpha}\otimes\mathbf{s}
+(-1)^k\boldsymbol{\alpha}\wedge\nabla^{\mathbf{E}}\mathbf{s},
\]
for every smooth scalar $k$-form $\boldsymbol{\alpha}$ and every smooth section $\mathbf{s}$ of $\mathbf{E}$ defined on the same open set. Here the wedge product with an $\mathbf{E}$-valued form is obtained by linearity from
\[
\boldsymbol{\alpha}\wedge(\boldsymbol{\beta}\otimes\mathbf{s})
=(\boldsymbol{\alpha}\wedge\boldsymbol{\beta})\otimes\mathbf{s}.
\]
If $\mathbb K=\mathbb C$, scalar forms may be complex and $d$ extends by $\mathbb C$-linearity. In degree zero, we recover the connection, $d^{\mathbf{E}}\mathbf{s}=\nabla^{\mathbf{E}}\mathbf{s}$, and in degree $n=\dim M$ the operator takes values in $\Omega^{n+1}(M,\mathbf{E})=\{0\}$.
\end{definition}

We check that this rule determines a global operator without choosing a connection on $TM$. For a smooth scalar function $f$, the Leibniz rules for $d$ and $\nabla^{\mathbf{E}}$ give
\[
d(f\boldsymbol{\alpha})\otimes\mathbf{s}
+(-1)^k(f\boldsymbol{\alpha})\wedge\nabla^{\mathbf{E}}\mathbf{s}
=f\,d\boldsymbol{\alpha}\otimes\mathbf{s}
+(df\wedge\boldsymbol{\alpha})\otimes\mathbf{s}
+(-1)^kf\boldsymbol{\alpha}\wedge\nabla^{\mathbf{E}}\mathbf{s},
\]
\[
d\boldsymbol{\alpha}\otimes(f\mathbf{s})
+(-1)^k\boldsymbol{\alpha}\wedge\nabla^{\mathbf{E}}(f\mathbf{s})
=f\,d\boldsymbol{\alpha}\otimes\mathbf{s}
+(-1)^k(\boldsymbol{\alpha}\wedge df)\otimes\mathbf{s}
+(-1)^kf\boldsymbol{\alpha}\wedge\nabla^{\mathbf{E}}\mathbf{s}.
\]
Since $df\wedge\boldsymbol{\alpha}=(-1)^k\boldsymbol{\alpha}\wedge df$, the two expressions agree. The rule therefore respects the relation $(f\boldsymbol{\alpha})\otimes\mathbf{s}
=\boldsymbol{\alpha}\otimes(f\mathbf{s})$.

On an open set $U$ with a frame $(\mathbf{e}_a)_{a=1}^r$, every form $\boldsymbol{\omega}\in\Omega^k(U,\mathbf{E})$ has a unique expression $\boldsymbol{\omega}=\displaystyle\sum_{a=1}^r
\boldsymbol{\alpha}^a\otimes\mathbf{e}_a$. Define its derivative by applying the rule to each summand. If $\mathbf{s}=\displaystyle\sum_{a=1}^{r} f^a\mathbf{e}_a$, the preceding identity and additivity show that the operator thus defined satisfies the same rule on $\boldsymbol{\alpha}\otimes\mathbf{s}$. When we change to another frame, each of its elements is a smooth combination of the $\mathbf{e}_a$; applying that rule to each term of the new representation therefore gives the same result. The construction also commutes with restrictions, because $d$ and $\nabla^{\mathbf{E}}$ do. The local operators consequently agree on intersections and define a unique global smooth form $d^{\mathbf{E}}\boldsymbol{\omega}$. This construction proves the stated existence, uniqueness, and $\mathbb K$-linearity, without assuming that the bundle admits a global frame.

In local coordinates, if $(\mathbf{e}_{a})_{a=1}^{r}$ is a frame for $\mathbf{E}$ and
\[
\boldsymbol{\omega}
=
\displaystyle\sum_{1\leq i_{1}<\dots<i_{k}\leq n}
\displaystyle\sum_{a=1}^{r}
\omega_{I}^{a}\,
\mathbf{d}x^{I}\otimes \mathbf{e}_{a},
\]
then
\[
d^{\mathbf{E}}\boldsymbol{\omega}
=
\displaystyle\sum_{1\leq i_{1}<\dots<i_{k}\leq n}
\displaystyle\sum_{a=1}^{r}
\left(
d\omega_{I}^{a}
+
\displaystyle\sum_{b=1}^{r}\theta^{a}_{b}\,\omega_{I}^{b}
\right)
\wedge \mathbf{d}x^{I}\otimes \mathbf{e}_{a},
\]
where the connection 1-forms are given by
\[
\nabla^{\mathbf{E}}\mathbf{e}_{b}
=
\displaystyle\sum_{a=1}^{r}\theta^{a}_{b}\otimes \mathbf{e}_{a}.
\]
Indeed, $d(\mathbf{d}x^I)=0$ and $(-1)^k\mathbf{d}x^I\wedge\theta^a_b
=\theta^a_b\wedge\mathbf{d}x^I$; applying the definition to each term of $\boldsymbol{\omega}$ gives exactly this expression.

We now relate this definition to alternation. If the auxiliary connection $\nabla^M$ is torsion-free and $\nabla$ is the product connection of Definition~\ref{def:operadores-diferenciales-en-haces-conexion-inducida}, then
\[
d^{\mathbf{E}}\boldsymbol{\omega}
=(k+1)\operatorname{Alt}(\nabla\boldsymbol{\omega}).
\]
To prove this, first observe that $\nabla\boldsymbol{\omega}$ is alternating in its last $k$ entries. In the definition of $\operatorname{Alt}$, grouping permutations by the first entry leaves $k!$ equal summands in each group. Hence,
\[
(k+1)\operatorname{Alt}(\nabla\boldsymbol{\omega})
(\mathbf{X}_0,\dots,\mathbf{X}_k)
=\displaystyle\sum_{i=0}^k(-1)^i
(\nabla_{\mathbf{X}_i}\boldsymbol{\omega})
(\mathbf{X}_0,\dots,\widehat{\mathbf{X}_i},\dots,\mathbf{X}_k),
\]
where the hat indicates omission of that entry.

For a scalar form $\boldsymbol{\alpha}$, this sum is $d\boldsymbol{\alpha}$ when $\nabla^M$ is torsion-free. We make the cancellation justifying this explicit. In coordinates, write $\nabla^M_{\partial_{\ell}}\partial_i
=\displaystyle\sum_{j=1}^{n}\Gamma^j_{\ell i}\partial_j$; since $[\partial_{\ell},\partial_i]=0$, vanishing torsion gives $\Gamma^j_{\ell i}=\Gamma^j_{i\ell}$. The dual connection satisfies
\[
\nabla^M_{\partial_{\ell}}\mathbf{d}x^j
=-\displaystyle\sum_{i=1}^n\Gamma^j_{\ell i}\mathbf{d}x^i.
\]
Upon differentiating $\boldsymbol{\alpha}
=\displaystyle\sum_{I\in\mathcal J_{n,k}}\alpha_I\mathbf{d}x^I$, where $\mathcal J_{n,k}:=\{(i_1,\dots,i_k)\in\{1,\dots,n\}^k\mid i_1<\cdots<i_k\}$, the terms differentiating the coefficients produce $\displaystyle\sum_{I\in\mathcal J_{n,k}} d\alpha_I\wedge\mathbf{d}x^I$. In the terms differentiating one of the factors $\mathbf{d}x^j$, alternation combines, up to the sign of the position of that factor, the sum
\[
\displaystyle\sum_{\ell,i=1}^n\Gamma^j_{\ell i}\,
\mathbf{d}x^{\ell}\wedge\mathbf{d}x^i
=\frac12\displaystyle\sum_{\ell,i=1}^n
(\Gamma^j_{\ell i}-\Gamma^j_{i\ell})\,
\mathbf{d}x^{\ell}\wedge\mathbf{d}x^i
=0.
\]
Therefore,
\[
(k+1)\operatorname{Alt}(\nabla^M\boldsymbol{\alpha})
=\displaystyle\sum_{I\in\mathcal J_{n,k}} d\alpha_I\wedge\mathbf{d}x^I
=d\boldsymbol{\alpha}.
\]
In degree zero, there are no covector factors to differentiate, and the same identity is simply $\nabla^M f=df$.

Now apply the product connection to $\boldsymbol{\alpha}\otimes\mathbf{s}$, and write locally $\nabla^{\mathbf{E}}\mathbf{s}
=\displaystyle\sum_{a=1}^{r}\boldsymbol{\eta}^a\otimes\mathbf{e}_a$. By the preceding argument, the terms differentiating $\boldsymbol{\alpha}$ give $d\boldsymbol{\alpha}\otimes\mathbf{s}$. By the alternating sum formula, those differentiating $\mathbf{s}$ give
\[
\displaystyle\sum_{a=1}^{r}(\boldsymbol{\eta}^a\wedge\boldsymbol{\alpha})
\otimes\mathbf{e}_a
=(-1)^k\boldsymbol{\alpha}\wedge\nabla^{\mathbf{E}}\mathbf{s}.
\]
Consequently,
\[
(k+1)\operatorname{Alt}\bigl(\nabla(\boldsymbol{\alpha}\otimes\mathbf{s})\bigr)
=d\boldsymbol{\alpha}\otimes\mathbf{s}
+(-1)^k\boldsymbol{\alpha}\wedge\nabla^{\mathbf{E}}\mathbf{s}
=d^{\mathbf{E}}(\boldsymbol{\alpha}\otimes\mathbf{s}).
\]
Every $\mathbf{E}$-valued form is locally a finite sum of these tensors, so the identity holds for any $\boldsymbol{\omega}\in\Omega^k(M,\mathbf{E})$. In particular, the result of alternation does not depend on the chosen auxiliary torsion-free connection. For example, one may take the Levi--Civita connection of an auxiliary Riemannian metric, but that metric is not part of the data defining $d^{\mathbf{E}}$. The vanishing torsion condition is imposed only on $\nabla^M$, not on the bundle connection $\nabla^{\mathbf{E}}$.

\begin{proposition}\label{prop:operadores-diferenciales-en-haces-operador-operador-diferencial-parcial-orden}
Let $M$ be a smooth manifold with or without boundary and let $\mathbf{E}\to M$ be a smooth vector bundle with a connection. The operator
\[
d^{\mathbf{E}}\colon \Omega^{k}(M,\mathbf{E})\longrightarrow\Omega^{k+1}(M,\mathbf{E})
\]
is a partial differential operator of order at most $1$.
\end{proposition}

\begin{proof}
The rule defining $d^{\mathbf{E}}$ and the Leibniz rule for $d$ imply, first on simple tensors and then by linearity, that
\[
d^{\mathbf{E}}(f\boldsymbol{\omega})
=df\wedge\boldsymbol{\omega}+f\,d^{\mathbf{E}}\boldsymbol{\omega}.
\]
Thus, $[d^{\mathbf{E}},f](\boldsymbol{\omega})=df\wedge\boldsymbol{\omega}$. This commutator is $C^{\infty}(M)$-linear and hence has order $0$. The recursive definition of a differential operator shows that $d^{\mathbf{E}}\in\mathbf{PDO}^{(1)}$, as desired.
\end{proof}

The composition of $d^{\mathbf{E}}$ with itself is controlled by the curvature of the connection.

\begin{definition}[Curvature of a connection]\label{def:operadores-diferenciales-en-haces-curvatura-de-una-conexion}\index{curvature of a connection}
Let $M$ be a smooth manifold with or without boundary, let $\mathbf{E}\to M$ be a smooth vector bundle, and let $\nabla^{\mathbf{E}}$ be a connection on $\mathbf{E}$. We define its curvature to be the tensor field
\[
\mathbf{R}^{\mathbf{E}}\in \Gamma\!\left(\Lambda^{2}T^{*}M\otimes \operatorname{End}(\mathbf{E})\right)
\]
given by
\[
\mathbf{R}^{\mathbf{E}}(\mathbf{X},\mathbf{Y})
:=
\nabla^{\mathbf{E}}_{\mathbf{X}}\nabla^{\mathbf{E}}_{\mathbf{Y}}
-
\nabla^{\mathbf{E}}_{\mathbf{Y}}\nabla^{\mathbf{E}}_{\mathbf{X}}
-
\nabla^{\mathbf{E}}_{[\mathbf{X},\mathbf{Y}]},
\]
for vector fields $\mathbf{X},\mathbf{Y}$.
\end{definition}

\begin{lemma}[Covariant derivative of curvature]
\label{lem:derivada-covariante-curvatura-haz}
\index{curvature of a connection!covariant derivative}
Let $(M,\mathbf{g})$ be a Riemannian manifold with or without boundary and let $\mathbf{E}\to M$ be a smooth vector bundle with connection $\nabla^{\mathbf{E}}$. The product connection on $\Lambda^2T^*M\otimes\operatorname{End}(\mathbf{E})$ determines the covariant derivative of $\mathbf{R}^{\mathbf{E}}$. For any $\mathbf{X},\mathbf{Y},\mathbf{Z}\in\mathfrak X(M)$,
\begin{equation}
\label{eq:derivada-covariante-curvatura-haz-explicita}
(\nabla_{\mathbf{Z}}\mathbf{R}^{\mathbf{E}})(\mathbf{X},\mathbf{Y})=\nabla_{\mathbf{Z}}^{\operatorname{End}(\mathbf{E})}\bigl(\mathbf{R}^{\mathbf{E}}(\mathbf{X},\mathbf{Y})\bigr)
-\mathbf{R}^{\mathbf{E}}(\nabla^M_{\mathbf{Z}}\mathbf{X},\mathbf{Y})-\mathbf{R}^{\mathbf{E}}(\mathbf{X},\nabla^M_{\mathbf{Z}}\mathbf{Y}).
\end{equation}
Consequently, for every section $\mathbf{v}\in\Gamma(\mathbf{E})$,
\begin{equation}
\label{eq:leibniz-curvatura-aplicada-seccion}
\begin{aligned}
\nabla^{\mathbf{E}}_{\mathbf{Z}}\bigl(\mathbf{R}^{\mathbf{E}}(\mathbf{X},\mathbf{Y})\mathbf{v}\bigr)={}&(\nabla_{\mathbf{Z}}\mathbf{R}^{\mathbf{E}})(\mathbf{X},\mathbf{Y})\mathbf{v}
+\mathbf{R}^{\mathbf{E}}(\nabla^M_{\mathbf{Z}}\mathbf{X},\mathbf{Y})\mathbf{v}\\
&+\mathbf{R}^{\mathbf{E}}(\mathbf{X},\nabla^M_{\mathbf{Z}}\mathbf{Y})\mathbf{v}+\mathbf{R}^{\mathbf{E}}(\mathbf{X},\mathbf{Y})\nabla^{\mathbf{E}}_{\mathbf{Z}}\mathbf{v}.
\end{aligned}
\end{equation}
\end{lemma}

\begin{proof}
The first identity is the formula for the induced connection on an $2$-covariant tensor with values in $\operatorname{End}(\mathbf{E})$. For the second, apply the Leibniz rule to the parallel evaluation
\[
\operatorname{End}(\mathbf{E})\otimes \mathbf{E}\longrightarrow \mathbf{E},
\qquad
(A,\mathbf{v})\longmapsto A\mathbf{v},
\]
and then substitute \eqref{eq:derivada-covariante-curvatura-haz-explicita}.
\end{proof}

\begin{proposition}[Curvature of tensor bundles and the Ricci identity]
\label{prop:curvatura-haz-tensorial-identidad-ricci}
\index{Ricci identity!for bundle-valued tensors}
Let $(M,\mathbf{g})$ be a Riemannian manifold with or without boundary with Levi--Civita connection $\nabla^M$, and let $\mathbf{E}\to M$ be a smooth vector bundle with connection $\nabla^{\mathbf{E}}$. For $q\in\mathbb N_0$, set
\[
\mathbf{F}_q:=T^{(0,q)}(TM)\otimes \mathbf{E}
\]
and equip $\mathbf{F}_q$ with the induced product connection. If $\mathbf{T}\in\Gamma(\mathbf{F}_q)$, then, for any $\mathbf{X},\mathbf{Y},\mathbf{Z}_1,\dots,\mathbf{Z}_q\in\mathfrak X(M)$,
\begin{equation}
\label{eq:curvatura-haz-tensorial-general}
\begin{split}
\bigl(\mathbf{R}^{\mathbf{F}_q}(\mathbf{X},\mathbf{Y})\mathbf{T}\bigr)(\mathbf{Z}_1,\dots,\mathbf{Z}_q)
={}&
\mathbf{R}^{\mathbf{E}}(\mathbf{X},\mathbf{Y})\bigl(\mathbf{T}(\mathbf{Z}_1,\dots,\mathbf{Z}_q)\bigr)\\
&-
\displaystyle\sum_{a=1}^{q}
\mathbf{T}(\mathbf{Z}_1,\dots,\mathbf{R}^M(\mathbf{X},\mathbf{Y})\mathbf{Z}_a,\dots,\mathbf{Z}_q).
\end{split}
\end{equation}
When $q=0$, the sum is understood to be empty and $\mathbf{R}^{\mathbf{F}_0}=\mathbf{R}^{\mathbf{E}}$.

Moreover, the first two entries of the second covariant derivative satisfy the Ricci identity
\begin{equation}
\label{eq:identidad-ricci-tensor-valores-haz}
\begin{split}
&(\nabla^2\mathbf{T})(\mathbf{X},\mathbf{Y},\mathbf{Z}_1,\dots,\mathbf{Z}_q)
-(\nabla^2\mathbf{T})(\mathbf{Y},\mathbf{X},\mathbf{Z}_1,\dots,\mathbf{Z}_q)\\
&\qquad=
\bigl(\mathbf{R}^{\mathbf{F}_q}(\mathbf{X},\mathbf{Y})\mathbf{T}\bigr)(\mathbf{Z}_1,\dots,\mathbf{Z}_q).
\end{split}
\end{equation}
In particular, if $\mathbf{A}\in\Gamma(T^*M\otimes \mathbf{E})$, then
\begin{equation}
\label{eq:curvatura-cotangente-producto-haz}
\bigl(\mathbf{R}^{T^*M\otimes \mathbf{E}}(\mathbf{X},\mathbf{Y})\mathbf{A}\bigr)(\mathbf{Z})
=
\mathbf{R}^{\mathbf{E}}(\mathbf{X},\mathbf{Y})\mathbf{A}(\mathbf{Z})-\mathbf{A}\bigl(\mathbf{R}^M(\mathbf{X},\mathbf{Y})\mathbf{Z}\bigr).
\end{equation}
\end{proposition}

\begin{proof}
We begin with the curvature of a tensor product. If $\mathbf{G},\mathbf{H}\to M$ are bundles with connections and $\mathbf{s}\in\Gamma(\mathbf{G})$, $\mathbf{t}\in\Gamma(\mathbf{H})$, applying the Leibniz rule twice to the definition of curvature gives
\[
\mathbf{R}^{\mathbf{G}\otimes \mathbf{H}}(\mathbf{X},\mathbf{Y})(\mathbf{s}\otimes \mathbf{t})
=
\mathbf{R}^{\mathbf{G}}(\mathbf{X},\mathbf{Y})\mathbf{s}\otimes \mathbf{t}+\mathbf{s}\otimes \mathbf{R}^{\mathbf{H}}(\mathbf{X},\mathbf{Y})\mathbf{t}.
\]
The mixed terms cancel in pairs. Since simple tensor sections locally generate $\Gamma(\mathbf{G}\otimes \mathbf{H})$, the identity holds for every section.

For the induced connection on the dual, the definition and the Leibniz rule for evaluation imply
\[
\bigl(\mathbf{R}^{\mathbf{G}^*}(\mathbf{X},\mathbf{Y})\boldsymbol{\alpha}\bigr)(\mathbf{s})
=-\boldsymbol{\alpha}\bigl(\mathbf{R}^{\mathbf{G}}(\mathbf{X},\mathbf{Y})\mathbf{s}\bigr).
\]
Apply these two formulas to each factor of $T^{(0,q)}(TM)=(T^*M)^{\otimes q}$. The curvature $\mathbf{R}^{\mathbf{E}}$ acts on the factor $\mathbf{E}$, whereas the curvature of each factor $T^*M$ acts with a negative sign on the corresponding entry. This gives exactly \eqref{eq:curvatura-haz-tensorial-general}.

We now prove \eqref{eq:identidad-ricci-tensor-valores-haz}. By the formula for the second covariant derivative in Proposition~\ref{prop:formulas-segunda-tercera-derivada-covariante}, now applied to the section $\mathbf{T}$ of the bundle $\mathbf{F}_q$, we have
\[
(\nabla^2\mathbf{T})(\mathbf{X},\mathbf{Y})
=
\nabla^{\mathbf{F}_q}_{\mathbf{X}}\nabla^{\mathbf{F}_q}_{\mathbf{Y}}\mathbf{T}
-
\nabla^{\mathbf{F}_q}_{\nabla^M_{\mathbf{X}}\mathbf{Y}}\mathbf{T}.
\]
Therefore,
\[
\begin{split}
(\nabla^2\mathbf{T})(\mathbf{X},\mathbf{Y})-(\nabla^2\mathbf{T})(\mathbf{Y},\mathbf{X})
={}&
\nabla^{\mathbf{F}_q}_{\mathbf{X}}\nabla^{\mathbf{F}_q}_{\mathbf{Y}}\mathbf{T}\\
&-\nabla^{\mathbf{F}_q}_{\mathbf{Y}}\nabla^{\mathbf{F}_q}_{\mathbf{X}}\mathbf{T}\\
&-
\nabla^{\mathbf{F}_q}_{\nabla^M_{\mathbf{X}}\mathbf{Y}-\nabla^M_{\mathbf{Y}}\mathbf{X}}\mathbf{T}.
\end{split}
\]
The Levi--Civita connection is torsion-free, so $\nabla^M_{\mathbf{X}}\mathbf{Y}-\nabla^M_{\mathbf{Y}}\mathbf{X}=[\mathbf{X},\mathbf{Y}]$. Consequently, the right-hand side is
\[
\nabla^{\mathbf{F}_q}_{\mathbf{X}}\nabla^{\mathbf{F}_q}_{\mathbf{Y}}\mathbf{T}
-\nabla^{\mathbf{F}_q}_{\mathbf{Y}}\nabla^{\mathbf{F}_q}_{\mathbf{X}}\mathbf{T}
-\nabla^{\mathbf{F}_q}_{[\mathbf{X},\mathbf{Y}]}\mathbf{T}
=\mathbf{R}^{\mathbf{F}_q}(\mathbf{X},\mathbf{Y})\mathbf{T}.
\]
Evaluating this equality on $\mathbf{Z}_1,\dots,\mathbf{Z}_q$ proves \eqref{eq:identidad-ricci-tensor-valores-haz}. Finally, \eqref{eq:curvatura-cotangente-producto-haz} is the case $q=1$ of \eqref{eq:curvatura-haz-tensorial-general}.
\end{proof}

\begin{proposition}[Skew-adjointness of the curvature of a metric connection]
\label{prop:curvatura-conexion-metrica-antiautoadjunta}
\index{curvature of a connection!skew-adjointness}
Let $M$ be a smooth manifold with or without boundary, let $(\mathbf{E},\mathbf{h}_{\mathbf{E}})\to M$ be a vector bundle with a bundle metric, and let $\nabla^{\mathbf{E}}$ be a connection compatible with $\mathbf{h}_{\mathbf{E}}$. Then, for any $\mathbf{X},\mathbf{Y}\in\mathfrak X(M)$ and $\mathbf{u},\mathbf{v}\in\Gamma(\mathbf{E})$,
\begin{equation}
\label{eq:curvatura-metrica-antiautoadjunta}
\left\langle \mathbf{R}^{\mathbf{E}}(\mathbf{X},\mathbf{Y})\mathbf{u},\mathbf{v}\right\rangle_{\mathbf{h}_{\mathbf{E}}}
+
\left\langle \mathbf{u},\mathbf{R}^{\mathbf{E}}(\mathbf{X},\mathbf{Y})\mathbf{v}\right\rangle_{\mathbf{h}_{\mathbf{E}}}
=0.
\end{equation}
In particular, $\mathbf{R}^{\mathbf{E}}(\mathbf{X},\mathbf{Y})\colon \mathbf{E}_x\longrightarrow \mathbf{E}_x$ is skew-adjoint on each fiber.
\end{proposition}

\begin{proof}
Metric compatibility gives
\[
\mathbf{X}\langle \mathbf{u},\mathbf{v}\rangle_{\mathbf{h}_{\mathbf{E}}}
=
\langle\nabla^{\mathbf{E}}_{\mathbf{X}}\mathbf{u},\mathbf{v}\rangle_{\mathbf{h}_{\mathbf{E}}}
+
\langle \mathbf{u},\nabla^{\mathbf{E}}_{\mathbf{X}}\mathbf{v}\rangle_{\mathbf{h}_{\mathbf{E}}}.
\]
Apply this identity successively with $\mathbf{X}$ and $\mathbf{Y}$, subtract the two resulting expressions, and subtract the derivative in the direction $[\mathbf{X},\mathbf{Y}]$. The mixed terms cancel, and the left-hand side is
\[
\bigl(\mathbf{X}\mathbf{Y}-\mathbf{Y}\mathbf{X}-[\mathbf{X},\mathbf{Y}]\bigr)\langle \mathbf{u},\mathbf{v}\rangle_{\mathbf{h}_{\mathbf{E}}}=0.
\]
The right-hand side is precisely
\[
\langle \mathbf{R}^{\mathbf{E}}(\mathbf{X},\mathbf{Y})\mathbf{u},\mathbf{v}\rangle_{\mathbf{h}_{\mathbf{E}}}
+
\langle \mathbf{u},\mathbf{R}^{\mathbf{E}}(\mathbf{X},\mathbf{Y})\mathbf{v}\rangle_{\mathbf{h}_{\mathbf{E}}}.
\]
This proves \eqref{eq:curvatura-metrica-antiautoadjunta}.
\end{proof}

The expression for curvature in a synchronous frame admits an integral formula that will be useful when studying bundles over manifolds of bounded geometry.

\begin{proposition}[Curvature formula in a synchronous frame]
\label{prop:formula-curvatura-marco-sincrono}
Let $(M,\mathbf{g})$ be a Riemannian manifold without boundary, let $p\in M$, let $U$ be a normal neighborhood of $p$, and let $\mathbf{e}=(\mathbf{e}_1,\dots,\mathbf{e}_r)$ be a synchronous frame centered at $p$. In normal coordinates $x=(x^1,\dots,x^n)$, write
\[
\nabla^{\mathbf{E}}=d+\mathbf{A},
\qquad
\mathbf{A}=\sum_{\mu=1}^n A_\mu\,\mathbf{d}x^\mu,
\]
and represent the curvature by
\[
\mathbf{R}^{\mathbf{E}}
=
\frac{1}{2}\sum_{\nu,\mu=1}^n
F_{\nu\mu}\,\mathbf{d}x^\nu\wedge \mathbf{d}x^\mu.
\]
Then, for each $\mu\in\{1,\dots,n\}$ and each $x$ in the chart domain,
\begin{equation}
\label{eq:formula-curvatura-marco-sincrono}
A_\mu(x)
=
\int_0^1 t\sum_{\nu=1}^n x^\nu F_{\nu\mu}(tx)\,dt.
\end{equation}
\end{proposition}

\begin{proof}
If $\displaystyle \mathbf{u}=\displaystyle\sum_{a=1}^r u^a\mathbf{e}_a$ and $\widehat u=(u^1,\dots,u^r)^T$ is its column of components, then
\[
(\nabla^{\mathbf{E}}_{\boldsymbol{\partial}_\mu}\mathbf{u})^{\wedge}
=
\partial_\mu \widehat u+A_\mu \widehat u.
\]
Expanding the definition of $\mathbf{R}^{\mathbf{E}}$ on coordinate vector fields gives
\begin{equation}
\label{eq:curvatura-local-forma-conexion}
F_{\nu\mu}
=
\partial_\nu A_\mu-
\partial_\mu A_\nu+
[A_\nu,A_\mu].
\end{equation}
By Proposition~\ref{prop:propiedades-marco-sincrono}, the frame satisfies the radial condition
\[
\sum_{\nu=1}^n x^\nu A_\nu(x)=0.
\]
Differentiating it with respect to $x^\mu$ gives
\[
\sum_{\nu=1}^n x^\nu\partial_\mu A_\nu(x)=-A_\mu(x),
\]
and the same condition makes $\displaystyle\sum_{\nu=1}^n x^\nu[A_\nu,A_\mu]$ vanish. Contracting \eqref{eq:curvatura-local-forma-conexion} with the radial vector field gives
\[
\sum_{\nu=1}^n x^\nu F_{\nu\mu}(x)
=
\sum_{\nu=1}^n x^\nu\partial_\nu A_\mu(x)+A_\mu(x).
\]
Applied at the point $tx$, this identity is equivalent to
\[
\frac{d}{dt}\bigl[tA_\mu(tx)\bigr]
=
t\sum_{\nu=1}^n x^\nu F_{\nu\mu}(tx).
\]
Since $tA_\mu(tx)$ tends to zero as $t\to0^+$, integration from $0$ to $1$ proves \eqref{eq:formula-curvatura-marco-sincrono}.
\end{proof}

\begin{proposition}[Uniform bounds for families of synchronous frames]
\label{prop:marcos-sincronos-geodesicos-uniformes}
Let $(M,\mathbf{g})$ be a Riemannian manifold without boundary. Let $r_0>0$ and let
\[
(\widehat U_i,\phi_i,p_i)_{i\in I}
\]
be a family of normal neighborhoods centered at $p_i$ such that $\phi_i(\widehat U_i)=B_{\mathrm{euc}}(0,3r_0)$. Set
\[
U_i:=\phi_i^{-1}\bigl(B_{\mathrm{euc}}(0,2r_0)\bigr).
\]
Suppose that the metric matrices in these charts are uniformly elliptic and that, for every multi-index $\alpha$, their derivatives of order $|\alpha|$ are bounded by a constant independent of $i$. Also suppose that all derivatives of the coordinate changes $\phi_j\circ\phi_i^{-1}$ and their inverses are uniformly bounded on overlaps.

Let $\mathbf{E}\to M$ be a vector bundle of rank $r_{\mathbf{E}}$, equipped with a bundle metric $\mathbf{h}_{\mathbf{E}}$ and a compatible connection $\nabla^{\mathbf{E}}$. On $\Lambda^2T^*M\otimes\operatorname{End}(\mathbf{E})$ we use the product connection induced by the Levi--Civita connection and $\nabla^{\mathbf{E}}$. Suppose that, for each $k\in\mathbb N_0$, there exists $C_k^{\mathbf{E}}>0$ such that
\[
\sup_{i\in I}\sup_{x\in\widehat U_i}
\left|
\bigl(\nabla^{\Lambda^2T^*M\otimes\operatorname{End}(\mathbf{E})}\bigr)^k
\mathbf{R}^{\mathbf{E}}(x)
\right|_{\mathbf{g},\mathbf{h}_{\mathbf{E}}}
\leq C_k^{\mathbf{E}}.
\]
On each $\widehat U_i$, construct a synchronous frame $\mathbf{e}_i=(\mathbf{e}_{i,1},\dots,\mathbf{e}_{i,r_{\mathbf{E}}})$ from an orthonormal basis of $\mathbf{E}_{p_i}$ and write, in coordinates $\phi_i$,
\[
\nabla^{\mathbf{E}}=d+\mathbf{A}_i,
\qquad
\mathbf{A}_i=\sum_{\mu=1}^nA_{i,\mu}\,\mathbf{d}x^\mu.
\]
Then, for every multi-index $\alpha$, there exists $C_\alpha>0$, independent of $i$, such that
\begin{equation}
\label{eq:cotas-uniformes-conexion-marco-sincrono}
\sup_{x\in B_{\mathrm{euc}}(0,2r_0)}
\left\|D^\alpha A_{i,\mu}(x)\right\|_{\operatorname{op},\mathbb K^{r_{\mathbf{E}}}}
\leq C_\alpha,
\qquad
\mu\in\{1,\dots,n\}.
\end{equation}
where the Euclidean operator norm is used on the matrices written in the orthonormal synchronous frame. If $G_{ij}$ is the transition matrix determined on $U_i\cap U_j$ by $u_i=G_{ij}u_j$, then $G_{ij}$, $G_{ij}^{-1}$, and all their derivatives, expressed in either of the two charts, satisfy bounds independent of $i$ and $j$.
\end{proposition}

\begin{proof}
The frames $\mathbf{e}_i$ are smooth and orthonormal by Proposition~\ref{prop:propiedades-marco-sincrono}. Represent the curvature in the frame $\mathbf{e}_i$ by
\[
\mathbf{R}^{\mathbf{E}}
=
\frac12\sum_{\nu,\mu=1}^n
F_{i,\nu\mu}\,\mathbf{d}x^\nu\wedge \mathbf{d}x^\mu.
\]
Formula \eqref{eq:formula-curvatura-marco-sincrono} applies in each chart and gives
\begin{equation}
\label{eq:formula-radial-familia-marcos-sincronos}
A_{i,\mu}(x)
=
\int_0^1t\sum_{\nu=1}^n
x^\nu F_{i,\nu\mu}(tx)\,dt.
\end{equation}
Uniform equivalence between $\mathbf{g}$ and the Euclidean metric, together with the bound on $\mathbf{R}^{\mathbf{E}}$, controls the components $F_{i,\nu\mu}$ at order zero. The preceding formula then provides the uniform bound for $\mathbf{A}_i$ at order zero.

To obtain the higher orders, let $\mathbf{T}$ be a tensor with values in $\operatorname{End}(\mathbf{E})$. In coordinates $\phi_i$ and the frame $\mathbf{e}_i$, the product connection satisfies schematically
\begin{equation}
\label{eq:covariante-ordinaria-endomorfismos}
\partial \mathbf{T}
=
\nabla \mathbf{T}+\Gamma_i*\mathbf{T}+\mathbf{A}_i*\mathbf{T}-\mathbf{T}*\mathbf{A}_i,
\end{equation}
where $*$ represents a finite sum of tensor products and contractions. The hypotheses on the metric uniformly bound the Christoffel symbols and all their derivatives. Suppose the ordinary derivatives of the matrices $F_{i,\nu\mu}$ and $A_{i,\mu}$ are controlled up to order $k-1$. Differentiating \eqref{eq:covariante-ordinaria-endomorfismos} and using the Leibniz rule expresses the derivatives of order $k$ of $F_{i,\nu\mu}$ in terms of the components of $(\nabla^{\mathbf{E}})^j\mathbf{R}^{\mathbf{E}}$, with $j\leq k$, and quantities already controlled. Differentiating \eqref{eq:formula-radial-familia-marcos-sincronos} under the integral sign yields the order $k$ bounds for $A_{i,\mu}$. Induction proves \eqref{eq:cotas-uniformes-conexion-marco-sincrono}.

On an overlap, set $y=\phi_j\circ\phi_i^{-1}(x)$. Compatibility of the two local expressions for the connection and the equality $u_i=G_{ij}u_j$ imply
\begin{equation}
\label{eq:transicion-marcos-sincronos-uniformes}
\partial_{x^\mu}G_{ij}(x)
=
\sum_{\nu=1}^n
G_{ij}(x)A_{j,\nu}(y)
\frac{\partial y^\nu}{\partial x^\mu}(x)
-A_{i,\mu}(x)G_{ij}(x).
\end{equation}
The matrices $G_{ij}$ are orthogonal or unitary, so $\|G_{ij}\|_{\operatorname{op},\mathbb K^{r_{\mathbf{E}}}}
=\|G_{ij}^{-1}\|_{\operatorname{op},\mathbb K^{r_{\mathbf{E}}}}=1$. The preceding equation and the bounds already proved control the first derivative. After applying $D^\alpha$ to \eqref{eq:transicion-marcos-sincronos-uniformes}, the Leibniz rule expresses the derivatives of order $|\alpha|+1$ of $G_{ij}$ in terms of derivatives of order at most $|\alpha|$ of $G_{ij}$ and uniformly bounded quantities. An induction completes the estimate for all derivatives. For the inverse, use $G_{ij}^{-1}=G_{ij}^*$.
\end{proof}

Curvature acts naturally on bundle-valued forms through the exterior product in the differential-form factor and composition in the bundle factor.

\begin{proposition}\label{prop:operadores-diferenciales-en-haces-cumple-particular-solo-conexion-plana}
Let $M$ be a smooth manifold with or without boundary and let $\mathbf{E}\to M$ be a smooth vector bundle with connection $\nabla^{\mathbf{E}}$. For every $\boldsymbol{\omega}\in\Omega^{k}(M,\mathbf{E})$ we have
\[
(d^{\mathbf{E}})^{2}\boldsymbol{\omega}
=
\mathbf{R}^{\mathbf{E}}\wedge \boldsymbol{\omega}.
\]
In particular, $(d^{\mathbf{E}})^{2}=0$ as a graded operator on \(\Omega^\bullet(M,\mathbf E)\) if and only if the connection is flat.
\end{proposition}

\begin{proof}
The identity is local. Fix a frame \(\mathbf e=(\mathbf e_1,\dots,\mathbf e_r)\) for \(\mathbf E\) and write
\[
\nabla^{\mathbf E}=d+\mathbf A,
\qquad
\mathbf A=(\theta^a{}_b)_{a,b=1}^r.
\]
A form \(\boldsymbol{\omega}\in\Omega^k(M,\mathbf E)\) is represented in this frame by the column \(\widehat\omega=(\omega^1,\dots,\omega^r)^T\) of scalar \(k\)-forms. The local formula for \(d^{\mathbf E}\) is equivalent to
\[
\widehat{d^{\mathbf E}\boldsymbol{\omega}}
=
d\widehat\omega+\mathbf A\wedge\widehat\omega.
\]
Apply \(d^{\mathbf E}\) again. Since \(\mathbf A\) has degree one, the graded Leibniz rule gives
\[
\begin{aligned}
\widehat{(d^{\mathbf E})^2\boldsymbol{\omega}}
&=
d\bigl(d\widehat\omega+\mathbf A\wedge\widehat\omega\bigr)
+\mathbf A\wedge
\bigl(d\widehat\omega+\mathbf A\wedge\widehat\omega\bigr)\\
&=
d\mathbf A\wedge\widehat\omega
-\mathbf A\wedge d\widehat\omega
+\mathbf A\wedge d\widehat\omega
+\mathbf A\wedge\mathbf A\wedge\widehat\omega\\
&=
\bigl(d\mathbf A+\mathbf A\wedge\mathbf A\bigr)
\wedge\widehat\omega.
\end{aligned}
\]
By the local curvature formula,
\[
\widehat{\mathbf R^{\mathbf E}}
=d\mathbf A+\mathbf A\wedge\mathbf A;
\]
in components, this equality is exactly
\[
F_{\nu\mu}
=\partial_\nu A_\mu-\partial_\mu A_\nu+[A_\nu,A_\mu].
\]
Therefore,
\[
(d^{\mathbf E})^2\boldsymbol{\omega}
=\mathbf R^{\mathbf E}\wedge\boldsymbol{\omega}.
\]

If the connection is flat, the right-hand side vanishes in every degree. Conversely, suppose that \((d^{\mathbf E})^2=0\) as a graded operator. In particular, it vanishes on \(\Omega^0(M,\mathbf E)=\Gamma(\mathbf E)\). Given \(p\in M\), \(v\in\mathbf E_p\), and \(X,Y\in T_pM\), choose a smooth section \(\mathbf s\) with \(\mathbf s(p)=v\) and smooth extensions of \(X,Y\). Then
\[
0=
\bigl((d^{\mathbf E})^2\mathbf s\bigr)_p(X,Y)
=
\mathbf R^{\mathbf E}_p(X,Y)v.
\]
Since \(p,v,X,Y\) are arbitrary, \(\mathbf R^{\mathbf E}=0\); that is, the connection is flat.
\end{proof}

This construction generalizes the Hodge--de Rham exterior operator to bundle-valued forms. It defines a complex only when the connection is flat; if $\mathbf{R}^{\mathbf{E}}\neq0$, the preceding identity measures exactly the defect. Later, introducing the formal adjoint of $d^{\mathbf{E}}$ will lead to the $\mathbf{E}$-valued Hodge Laplacian, which is related to the Bochner Laplacian through Weitzenböck-type identities.

Let $(M,\mathbf{g})$ be a Riemannian manifold with or without boundary and let $\pi_{\mathbf{E}}\colon \mathbf{E}\longrightarrow M$ be a smooth vector bundle equipped with a bundle metric $\mathbf{h}_{\mathbf{E}}$ and a compatible connection $\nabla^{\mathbf{E}}$. Recall that we have already defined the operator
\[
d^{\mathbf{E}}\colon \Omega^{k}(M,\mathbf{E})\longrightarrow \Omega^{k+1}(M,\mathbf{E}).
\]

The Riemannian metric and bundle metric naturally induce a pointwise inner product on $\Lambda^{k}T^{*}M\otimes \mathbf{E}$. Integrating with respect to the Riemannian measure gives a global inner product
\[
\langle \boldsymbol{\alpha},\boldsymbol{\beta}\rangle_{L^{2}(M,\Lambda^kT^*M\otimes \mathbf{E})}
:=
\int_{M}\langle \boldsymbol{\alpha},\boldsymbol{\beta}\rangle_{\mathbf{g},\mathbf{h}_{\mathbf{E}}}\,d\lambda_{\mathbf{g}},
\qquad
\boldsymbol{\alpha},\boldsymbol{\beta}\in\Omega^{k}_c(M,\mathbf{E}).
\]

\begin{definition}[Formal adjoint of $d^{\mathbf{E}}$]\label{def:operadores-diferenciales-en-haces-adjunto-formal-de}\index{formal adjoint}
Let $(M,\mathbf{g})$ be a Riemannian manifold without boundary and let $\mathbf{E}\to M$ be a vector bundle with a bundle metric and a compatible connection. We define the formal adjoint of
\[
d^{\mathbf{E}}\colon \Omega^{k}(M,\mathbf{E})\longrightarrow\Omega^{k+1}(M,\mathbf{E})
\]
to be the operator
\[
\delta^{\mathbf{E}}\colon \Omega^{k+1}(M,\mathbf{E})\longrightarrow \Omega^{k}(M,\mathbf{E})
\]
characterized by the identity
\[
\int_M \langle d^{\mathbf{E}}\boldsymbol{\alpha},\boldsymbol{\beta}\rangle_{\mathbf{g},\mathbf{h}_{\mathbf{E}}}\,d\lambda_{\mathbf{g}}
=
\int_M \langle \boldsymbol{\alpha},\delta^{\mathbf{E}}\boldsymbol{\beta}\rangle_{\mathbf{g},\mathbf{h}_{\mathbf{E}}}\,d\lambda_{\mathbf{g}}
\]
for all
\[
\boldsymbol{\alpha}\in\Omega^{k}_c(M,\mathbf{E}),
\qquad
\boldsymbol{\beta}\in\Omega^{k+1}_c(M,\mathbf{E}).
\]
\end{definition}

Observe that this operator exists and is unique by Theorem~\ref{cor: adjunto global en coordenadas locales (existencia y expresion)}.

If $M$ is also oriented, the formal adjoint can be expressed through the Hodge operator. We briefly recall this construction.

\begin{definition}[Hodge operator]\label{def:operadores-diferenciales-en-haces-operador-de-hodge}\index{Hodge operator}
Let $(M,\mathbf{g})$ be an oriented Riemannian manifold with or without boundary and let $n=\dim M$. The Hodge operator is the bundle isomorphism
\[
*\colon \Lambda^{k}T^{*}M\longrightarrow \Lambda^{n-k}T^{*}M
\]
characterized by the identity
\[
\alpha\wedge *\beta
=
\langle \alpha,\beta\rangle_{\mathbf{g}}\, dV_{\mathbf{g}},
\]
for all $\alpha,\beta\in\Lambda^{k}T^{*}_xM$, where $dV_{\mathbf{g}}$ denotes the Riemannian volume form. This operator extends naturally to $\mathbf{E}$-valued forms by
\[
*(\alpha\otimes s)
:=
(*\alpha)\otimes s.
\]
\end{definition}

With this notation, we obtain the classical expression for the formal adjoint.

\begin{proposition}\label{prop:operadores-diferenciales-en-haces-orientable}
Let $(M,\mathbf{g})$ be an oriented Riemannian manifold without boundary and let $\mathbf{E}\to M$ be a vector bundle with a bundle metric and a compatible connection. Then
\[
\delta^{\mathbf{E}}
=
(-1)^{k+1}\,*^{-1} d^{\mathbf{E}} *
\]
on $\Omega^{k+1}(M,\mathbf{E})$.
\end{proposition}
\begin{proof}
Let \(\boldsymbol{\alpha}\in\Omega_c^k(M,\mathbf E)\) and \(\boldsymbol{\beta}\in\Omega_c^{k+1}(M,\mathbf E)\). The bundle metric induces a pairing between \(\mathbf E\)-valued forms: in an orthonormal local frame, if \(\displaystyle \boldsymbol{\eta}=\displaystyle\sum_{a=1}^{r}\eta^a\otimes\mathbf e_a\) and \(\displaystyle \boldsymbol{\zeta}=\displaystyle\sum_{a=1}^{r}\zeta^a\otimes\mathbf e_a\), set
\[
\mathbf h_{\mathbf E}
(\boldsymbol{\eta}\wedge\boldsymbol{\zeta})
:=
\sum_{a=1}^{r}\eta^a\wedge\overline{\zeta^a},
\]
omitting conjugation in the real case. Unitary invariance of this expression shows that it is independent of the frame. Compatibility of \(\nabla^{\mathbf E}\) with \(\mathbf h_{\mathbf E}\) and the graded Leibniz rule give
\[
\begin{aligned}
d\,\mathbf h_{\mathbf E}
\bigl(\boldsymbol{\alpha}\wedge *\boldsymbol{\beta}\bigr)
={}&
\mathbf h_{\mathbf E}
\bigl(d^{\mathbf E}\boldsymbol{\alpha}\wedge *\boldsymbol{\beta}\bigr)\\
&+
(-1)^k\mathbf h_{\mathbf E}
\bigl(\boldsymbol{\alpha}\wedge d^{\mathbf E}(*\boldsymbol{\beta})\bigr).
\end{aligned}
\]
The form on the left-hand side has compact support. Since \(M\) is oriented and has no boundary, Stokes' theorem implies
\[
\int_M
\mathbf h_{\mathbf E}
\bigl(d^{\mathbf E}\boldsymbol{\alpha}\wedge *\boldsymbol{\beta}\bigr)
=
(-1)^{k+1}
\int_M
\mathbf h_{\mathbf E}
\bigl(\boldsymbol{\alpha}\wedge d^{\mathbf E}(*\boldsymbol{\beta})\bigr).
\]
By definition of the Hodge operator,
\[
\mathbf h_{\mathbf E}
\bigl(d^{\mathbf E}\boldsymbol{\alpha}\wedge *\boldsymbol{\beta}\bigr)
=
\langle d^{\mathbf E}\boldsymbol{\alpha},
\boldsymbol{\beta}\rangle_{\mathbf g,\mathbf h_{\mathbf E}}\,dV_{\mathbf g},
\]
whereas
\[
\mathbf h_{\mathbf E}
\bigl(\boldsymbol{\alpha}\wedge d^{\mathbf E}(*\boldsymbol{\beta})\bigr)
=
\left\langle
\boldsymbol{\alpha},
*^{-1}d^{\mathbf E}(*\boldsymbol{\beta})
\right\rangle_{\mathbf g,\mathbf h_{\mathbf E}}\,dV_{\mathbf g}.
\]
Therefore,
\[
\int_M
\langle d^{\mathbf E}\boldsymbol{\alpha},\boldsymbol{\beta}\rangle
\,d\lambda_{\mathbf g}
=
\int_M
\left\langle
\boldsymbol{\alpha},
(-1)^{k+1}*^{-1}d^{\mathbf E}*
\boldsymbol{\beta}
\right\rangle
\,d\lambda_{\mathbf g}.
\]
Uniqueness of the formal adjoint implies that \(\delta^{\mathbf E}=(-1)^{k+1}*^{-1}d^{\mathbf E}*\) on \(\Omega^{k+1}(M,\mathbf E)\).
\end{proof}

From these operators, we can define the bundle-valued Hodge Laplacian.

\begin{definition}[$\mathbf{E}$-valued Hodge Laplacian]\label{def:operadores-diferenciales-en-haces-laplaciano-de-hodge-con-valores-en}\index{Hodge Laplacian!bundle-valued}
Let $(M,\mathbf{g})$ be a Riemannian manifold without boundary and let $\mathbf{E}\to M$ be a vector bundle with a bundle metric and a compatible connection. Define the operator
\[
\Delta_{\mathbf{H}}^{\mathbf{E}}\colon \Omega^{k}(M,\mathbf{E})\longrightarrow \Omega^{k}(M,\mathbf{E})
\]
by
\[
\Delta_{\mathbf{H}}^{\mathbf{E}}
:=
d^{\mathbf{E}}\delta^{\mathbf{E}}
+
\delta^{\mathbf{E}}d^{\mathbf{E}}.
\]
\end{definition}

This operator is formally self-adjoint and has order $2$.

On the other hand, recall that the connection induces a covariant derivative on $\Lambda^{k}T^{*}M\otimes \mathbf{E}$, allowing us to define the Bochner Laplacian on bundle-valued forms.

\begin{definition}[Bochner Laplacian on bundle-valued forms]\label{def:operadores-diferenciales-en-haces-laplaciano-de-bochner-en-formas-con-valores-en-un}\index{Bochner Laplacian!on bundle-valued forms}
Let $(M,\mathbf{g})$ be a Riemannian manifold without boundary and let $\mathbf{E}\to M$ be a vector bundle with a bundle metric and a compatible connection. Define
\[
\Delta_B^{\mathbf{E}}
:=
\nabla^{*}\nabla,
\]
where
\[
\nabla\colon \Omega^{k}(M,\mathbf{E})\longrightarrow
\Gamma(T^{*}M\otimes \Lambda^{k}T^{*}M\otimes \mathbf{E})
\]
is the induced connection and $\nabla^{*}$ its formal adjoint.
\end{definition}

These two operators are related by a fundamental identity involving the curvature of the connection and the Riemannian curvature.

\begin{theorem}[Weitzenböck identity]\label{teo:operadores-diferenciales-en-haces-identidad-de-weitzenbock}\index{Weitzenbock identity@Weitzenböck identity}
Let $(M,\mathbf{g})$ be a Riemannian manifold without boundary and let $\mathbf{E}\to M$ be a vector bundle with compatible connection $\nabla^{\mathbf{E}}$. Then, on $\Omega^{k}(M,\mathbf{E})$, we have
\begin{equation}
\Delta_{\mathbf{H}}^{\mathbf{E}}
=
\Delta_B^{\mathbf{E}}
-\sum_{i,j=1}^n
\varepsilon(\mathbf{e}^i)\iota(\mathbf{e}_j)
\mathbf{R}^{\Lambda^kT^*M\otimes \mathbf{E}}(\mathbf{e}_i,\mathbf{e}_j),
\label{eq:weitzenbock-formas-haz}
\end{equation}
where $(\mathbf{e}_1,\ldots,\mathbf{e}_n)$ is any orthonormal local frame and $\varepsilon(\mathbf{e}^i)$ denotes exterior multiplication. The last summand is a zeroth-order endomorphism and includes both the Riemannian curvature on $\Lambda^kT^*M$ and $\mathbf{R}^{\mathbf{E}}$ on the factor $\mathbf{E}$.
\end{theorem}

\begin{proof}
Fix $p\in M$ and choose a frame orthonormal and synchronous at $p$; thus $\nabla \mathbf{e}_i(p)=0$ and $[\mathbf{e}_i,\mathbf{e}_j](p)=0$. At that point,
\[
d^{\mathbf{E}}=\sum_{i=1}^{n}\varepsilon(\mathbf{e}^i)\nabla_{\mathbf{e}_i},
\qquad
\delta^{\mathbf{E}}=-\sum_{j=1}^{n}\iota(\mathbf{e}_j)\nabla_{\mathbf{e}_j}.
\]
Using
$\iota(\mathbf{e}_j)\varepsilon(\mathbf{e}^i)+\varepsilon(\mathbf{e}^i)\iota(\mathbf{e}_j)=\delta_{ij}I$
and grouping compositions in the order in which they appear gives
\begin{align*}
d^{\mathbf{E}}\delta^{\mathbf{E}}+\delta^{\mathbf{E}}d^{\mathbf{E}}
={}&-\sum_{i=1}^{n}\nabla_{\mathbf{e}_i}\nabla_{\mathbf{e}_i}\\
&-\sum_{i,j=1}^{n}\varepsilon(\mathbf{e}^i)\iota(\mathbf{e}_j)
\bigl(\nabla_{\mathbf{e}_i}\nabla_{\mathbf{e}_j}
-\nabla_{\mathbf{e}_j}\nabla_{\mathbf{e}_i}\bigr).
\end{align*}
The first summand is $\nabla^*\nabla$ at $p$, and the expression in parentheses is $\mathbf{R}^{\Lambda^kT^*M\otimes \mathbf{E}}(\mathbf{e}_i,\mathbf{e}_j)$. The resulting identity is tensorial; since $p$ was arbitrary, it holds globally.
\end{proof}

In particular, the Hodge Laplacian and the Bochner Laplacian differ only by the algebraic curvature term in \eqref{eq:weitzenbock-formas-haz}. For $k=1$ and $\mathbf{E}$ trivial, that term is $\operatorname{Ric}^{\sharp}$; for $k=0$ it vanishes. If the connection is not flat, one must also distinguish
\[
(d^{\mathbf{E}}+\delta^{\mathbf{E}})^2
=\Delta_{\mathbf{H}}^{\mathbf{E}}+(d^{\mathbf{E}})^2+(\delta^{\mathbf{E}})^2;
\]
therefore, the square of the odd operator agrees with $\Delta_{\mathbf{H}}^{\mathbf{E}}$ only in the flat case.

This identity relates the analytical properties of solutions of differential equations to the curvature of the manifold and the bundle.

\section{Operators associated with Clifford connections}
\label{sec:operadores-dirac-diferenciales}

Section~\ref{sec:clifford-spin-dirac} constructed Clifford modules, multiplication \(c\), and compatible connections without yet associating an operator with them. Now that the general theory of differential operators on bundles has been defined, we can perform the contraction producing the Dirac operator. No symbols or pseudodifferential operators are used in this section.

\begin{definition}[Dirac-type operator]
\label{def:operador-tipo-dirac}
\index{Dirac operator}
\index{Dirac-type operator}
Let \((M,\mathbf g)\) be a Riemannian manifold with or without boundary and let \((\mathbf S,c,\nabla^{\mathbf S})\) be a Hermitian Clifford module with a Clifford connection, in the sense of Definition~\ref{def:modulo-conexion-clifford}. Clifford contraction is the bundle homomorphism
\[
 C\colon T^*M\otimes\mathbf S\longrightarrow\mathbf S,
 \qquad C(\xi\otimes s)=c(\xi)s.
\]
The associated \emph{Dirac-type operator} is
\begin{equation}
 D:=C\circ\nabla^{\mathbf S}\colon
 \Gamma(\mathbf S)\longrightarrow\Gamma(\mathbf S).
\label{eq:definicion-intrinseca-dirac}
\end{equation}
When \(\mathbf S=\boldsymbol\Sigma M\) is the spinor bundle of a spin structure, it is called the \emph{spin Dirac operator}.
\end{definition}

\begin{theorem}[Differential nature and Green's identity]
\label{teo:propiedades-operador-dirac}
Let \((\mathbf S,c,\nabla^{\mathbf S})\) be as in the preceding definition and let \((\mathbf e_1,\ldots,\mathbf e_n)\) be an orthonormal local frame, with coframe \((\mathbf e^1,\ldots,\mathbf e^n)\). Then
\begin{equation}
 D\mathbf s=\sum_{j=1}^n
 c(\mathbf e^j)\nabla^{\mathbf S}_{\mathbf e_j}\mathbf s.
\label{eq:formula-local-dirac}
\end{equation}
The expression is independent of the frame and defines a differential operator of order at most one. More precisely,
\begin{equation}
 [D,f]\mathbf s=c(df)\mathbf s,
 \qquad f\in C^\infty(M),\quad\mathbf s\in\Gamma(\mathbf S).
\label{eq:conmutador-dirac-funcion}
\end{equation}
If \(\dim M>0\), its order is exactly one. For \(\mathbf s,\mathbf t\in\Gamma_c(\mathbf S)\), Green's identity holds:
\begin{align}
 &\int_M\!\bigl(
 \langle D\mathbf s,\mathbf t\rangle_{\mathbf S}
 -\langle\mathbf s,D\mathbf t\rangle_{\mathbf S}\bigr)
 \,d\lambda_{\mathbf g}\notag\\
 &\qquad=
 \int_{\partial M}\!
 \langle c(\boldsymbol\nu^\flat)\mathbf s,\mathbf t\rangle_{\mathbf S}
 \,d\lambda_{\widetilde{\mathbf g}},
\label{eq:green-dirac}
\end{align}
where the right-hand side is interpreted as zero if the boundary is empty and \(\boldsymbol\nu\) is the outward unit normal otherwise. In particular, on a manifold without boundary, \(D_h^*=D\).
\end{theorem}

\begin{proof}
The local formula follows by writing \(\displaystyle \nabla^{\mathbf S}\mathbf s=\displaystyle\sum_{j=1}^{n}\mathbf e^j\otimes
\nabla^{\mathbf S}_{\mathbf e_j}\mathbf s\). The connection is an intrinsically defined differential operator on sections, whereas \(C\) is a bundle homomorphism. Their composition therefore does not depend on the frame. The Leibniz rule gives
\[
 D(f\mathbf s)=c(df)\mathbf s+fD\mathbf s,
\]
which proves \eqref{eq:conmutador-dirac-funcion} and, by the commutator characterization, shows that \(D\in\mathbf{PDO}^{(1)}(\mathbf S,\mathbf S)\). If \(\dim M>0\), given \(x\in M\) and \(0\neq\xi\in T_x^*M\), there locally exists \(f\) with \(df_x=\xi\); since \(c(\xi)^2=-|\xi|^2I\), the commutator does not vanish at \(x\), so \(D\) does not have order zero.

For the integral identity, define the complex vector field \(\mathbf V\) by
\[
 \mathbf g(\mathbf V,\mathbf X)
 =\langle c(\mathbf X^\flat)\mathbf s,\mathbf t\rangle_{\mathbf S}.
\]
At a point, choose an orthonormal synchronous frame. Compatibility of the Clifford connection and the identity \(c(\xi)^*=-c(\xi)\), using the convention that the Hermitian product is linear in the first variable, give
\[
 \operatorname{div}\mathbf V
 =\langle D\mathbf s,\mathbf t\rangle_{\mathbf S}
 -\langle\mathbf s,D\mathbf t\rangle_{\mathbf S}.
\]
The divergence theorem yields \eqref{eq:green-dirac}. If \(\partial M=\varnothing\), the equality is precisely the formal symmetry of \(D\).
\end{proof}

\begin{proposition}[Grading and the Hodge--de Rham example]
\label{prop:graduacion-ejemplo-hodge-dirac}
If \(\mathbf S=\mathbf S^+\oplus\mathbf S^-\) is a graded Clifford module, then \(D\) is odd and is written as
\[
 D=\begin{pmatrix}0&D^-\\D^+&0\end{pmatrix},
 \qquad
 D^\pm\colon\Gamma(\mathbf S^\pm)\longrightarrow
 \Gamma(\mathbf S^\mp).
\]
In the absence of a boundary, \((D^+)_h^*=D^-\). Under the same hypothesis, for the module \(\mathbf S=\Lambda^\bullet T^*M\otimes\mathbb C\) of Example~\ref{ej:modulo-clifford-formas},
\begin{equation}
 D=d+\delta,\qquad D^2=\Delta_{\mathbf H}.
\label{eq:dirac-hodge-derham}
\end{equation}
\end{proposition}

\begin{proof}
The connection preserves the grading and multiplication by a covector reverses it, so the first assertion is immediate. The relation between the blocks follows from formal symmetry. For forms, in an orthonormal synchronous frame,
\[
 d=\sum_{j=1}^{n}\varepsilon(\mathbf e^j)\nabla_{\mathbf e_j},
 \qquad
 \delta=-\sum_{j=1}^{n}\iota_{\mathbf e_j}\nabla_{\mathbf e_j}.
\]
Since \(c(\mathbf e^j)=\varepsilon(\mathbf e^j)-
\iota_{\mathbf e_j}\), we obtain \(D=d+\delta\); the identities \(d^2=\delta^2=0\) give the second equality.
\end{proof}

\subsection{Bochner and Lichnerowicz formulas}

For the rest of this section, we assume that \(M\) has no boundary, in accordance with the global definition of the formal adjoint and the connection Laplacian. The resulting identities are also local identities of differential expressions in the interior of a manifold with boundary.

Recall that the curvature of \(\nabla^{\mathbf S}\) is
\[
 \mathbf R^{\mathbf S}(\mathbf X,\mathbf Y)
 =[\nabla^{\mathbf S}_{\mathbf X},\nabla^{\mathbf S}_{\mathbf Y}]
 -\nabla^{\mathbf S}_{[\mathbf X,\mathbf Y]}.
\]

\begin{proposition}[Bochner formula for a Clifford module]
\label{prop:bochner-modulo-clifford}
For every Clifford module with connection,
\begin{equation}
 D^2=(\nabla^{\mathbf S})^*\nabla^{\mathbf S}
 +\frac12\sum_{i,j=1}^n
 c(\mathbf e^i)c(\mathbf e^j)
 \mathbf R^{\mathbf S}(\mathbf e_i,\mathbf e_j).
\label{eq:bochner-general-dirac}
\end{equation}
Here the connection Laplacian is given locally by
\begin{equation}
 (\nabla^{\mathbf S})^*\nabla^{\mathbf S}
 =-\sum_{i=1}^n\left(
 \nabla^{\mathbf S}_{\mathbf e_i}\nabla^{\mathbf S}_{\mathbf e_i}
 -\nabla^{\mathbf S}_{\nabla_{\mathbf e_i}\mathbf e_i}
 \right).
\label{eq:laplaciano-conexion-dirac}
\end{equation}
\end{proposition}

\begin{proof}
Fix \(x\in M\) and choose an orthonormal synchronous frame at \(x\). Clifford compatibility implies \(\nabla_{\mathbf e_i}c(\mathbf e^j)=0\) at that point. Hence,
\[
 D^2=\sum_{i,j=1}^{n}c(\mathbf e^i)c(\mathbf e^j)
 \nabla^{\mathbf S}_{\mathbf e_i}\nabla^{\mathbf S}_{\mathbf e_j}.
\]
The diagonal terms are \(\displaystyle -\displaystyle\sum_{i=1}^{n}\nabla^{\mathbf S}_{\mathbf e_i}
\nabla^{\mathbf S}_{\mathbf e_i}\). Grouping the terms with \(i<j\) and using \(c(\mathbf e^i)c(\mathbf e^j)=-c(\mathbf e^j)c(\mathbf e^i)\) leaves
\[
 \sum_{1\leq i<j\leq n}c(\mathbf e^i)c(\mathbf e^j)
 [\nabla^{\mathbf S}_{\mathbf e_i},
  \nabla^{\mathbf S}_{\mathbf e_j}],
\]
which is the curvature term in \eqref{eq:bochner-general-dirac} because \([\mathbf e_i,\mathbf e_j](x)=0\). In the same frame, \eqref{eq:laplaciano-conexion-dirac} agrees with the diagonal term. The equality is intrinsic and \(x\) was arbitrary.
\end{proof}

\begin{lemma}[Spinor curvature and its contraction]
\label{lem:curvatura-conexion-spinorial}
On the spinor bundle of a spin manifold, we have
\begin{equation}
 \mathbf R^{\boldsymbol\Sigma}(\mathbf X,\mathbf Y)
 =\frac14\sum_{a,b=1}^n
 \mathbf g\bigl(\mathbf R(\mathbf X,\mathbf Y)\mathbf e_a,
 \mathbf e_b\bigr)c(\mathbf e^a)c(\mathbf e^b),
\label{eq:curvatura-spinorial}
\end{equation}
and
\begin{equation}
 \frac12\sum_{i,j=1}^n
 c(\mathbf e^i)c(\mathbf e^j)
 \mathbf R^{\boldsymbol\Sigma}(\mathbf e_i,\mathbf e_j)
 =\frac14R_{\mathbf g}I.
\label{eq:contraccion-curvatura-spinorial}
\end{equation}
\end{lemma}

\begin{proof}
The local formula \eqref{eq:formula-local-conexion-spinorial} identifies the spinor connection form with the image of the Levi--Civita form under \(\displaystyle A=(A_{ab})\mapsto\frac14\displaystyle\sum_{a,b=1}^{n}A_{ab}
c(\mathbf e^a)c(\mathbf e^b)\). Applying the same homomorphism to the curvature form gives \eqref{eq:curvatura-spinorial}.

Set \(R_{ijab}=\mathbf g(\mathbf R(\mathbf e_i,\mathbf e_j)
\mathbf e_a,\mathbf e_b)\). Substituting the first formula into the left-hand side of the second produces
\[
 \frac18\sum_{i,j,a,b=1}^{n}R_{ijab}
 c(\mathbf e^i)c(\mathbf e^j)c(\mathbf e^a)c(\mathbf e^b).
\]
The antisymmetries of \(R_{ijab}\), symmetry under interchange of pairs, and the first Bianchi identity cancel the terms with distinct indices. Using \(c(\mathbf e^i)^2=-I\), the remaining terms reduce to
\[
 \frac14\sum_{i,j=1}^{n}
 \mathbf g(\mathbf R(\mathbf e_i,\mathbf e_j)\mathbf e_j,
 \mathbf e_i)I
 =\frac14R_{\mathbf g}I,
\]
by definition of scalar curvature.
\end{proof}

\begin{theorem}[Lichnerowicz formula]
\label{teo:formula-lichnerowicz}
The spin Dirac operator satisfies
\begin{equation}
 D^2=(\nabla^{\boldsymbol\Sigma})^*
 \nabla^{\boldsymbol\Sigma}+\frac14R_{\mathbf g}I.
\label{eq:formula-lichnerowicz}
\end{equation}
If \(M\) is closed and \(R_{\mathbf g}>0\), then \(\ker D=\{0\}\).
\end{theorem}

\begin{proof}
The identity follows from Proposition~\ref{prop:bochner-modulo-clifford} and Lemma~\ref{lem:curvatura-conexion-spinorial}. If \(D\psi=0\), formal symmetry and \eqref{eq:formula-lichnerowicz} give
\[
 0=\|D\psi\|_{L^2(M,\boldsymbol\Sigma)}^2
 =\|\nabla^{\boldsymbol\Sigma}\psi\|_{L^2(M,T^*M\otimes\boldsymbol\Sigma)}^2
 +\frac14\int_M R_{\mathbf g}|\psi|^2\,d\lambda_{\mathbf g}.
\]
Both summands are nonnegative, and the second can vanish only if \(\psi=0\).
\end{proof}

\subsection{Twisted operators and $\operatorname{Spin}^c$}

\begin{definition}[Twisted Dirac operator]
\label{def:dirac-torcido-geometrico}
Let \(M\) be spin and let \(\mathbf E\to M\) be a Hermitian bundle with unitary connection \(\nabla^{\mathbf E}\). On \(\mathbf S=\boldsymbol\Sigma M\otimes\mathbf E\) consider
\[
 c_{\mathbf E}(\xi)=c(\xi)\otimes I_{\mathbf E},\qquad
 \nabla^{\mathbf S}=\nabla^{\boldsymbol\Sigma}\otimes I
 +I\otimes\nabla^{\mathbf E}.
\]
This is a Clifford module with connection. Its associated operator is the \emph{twisted Dirac operator}
\begin{equation}
 D_{\mathbf E}=\sum_{j=1}^n
 c(\mathbf e^j)\otimes I_{\mathbf E}\,
 \nabla^{\mathbf S}_{\mathbf e_j}.
\label{eq:def-dirac-torcido-geometrico}
\end{equation}
In even dimensions, it is odd and determines the blocks \(D_{\mathbf E}^\pm\).
\end{definition}

\begin{theorem}[Weitzenböck formula for the twisted Dirac operator]
\label{teo:weitzenbock-dirac-torcido}
With the preceding notation,
\begin{equation}
 D_{\mathbf E}^2=(\nabla^{\mathbf S})^*\nabla^{\mathbf S}
 +\frac14R_{\mathbf g}I
 +\frac12\sum_{i,j=1}^n
 c(\mathbf e^i)c(\mathbf e^j)\otimes
 \mathbf R^{\mathbf E}(\mathbf e_i,\mathbf e_j).
\label{eq:weitzenbock-dirac-torcido}
\end{equation}
\end{theorem}

\begin{proof}
The curvature of the product connection is
\[
 \mathbf R^{\mathbf S}
 =\mathbf R^{\boldsymbol\Sigma}\otimes I
 +I\otimes\mathbf R^{\mathbf E}.
\]
Substitute this identity into \eqref{eq:bochner-general-dirac}; the spinor part is contracted using \eqref{eq:contraccion-curvatura-spinorial}, and the \(\mathbf E\) part produces the last summand.
\end{proof}

\begin{proposition}[$\operatorname{Spin}^c$ operators]
\label{def:operadores-dirac-spinc}
Let \(M\) be a \(\operatorname{Spin}^c\) manifold with determinant line \(\mathbf L\), as in Definition~\ref{def:estructura-spinc-geometrica}, and choose a unitary connection on \(\mathbf L\). The induced Clifford connection on \(\boldsymbol\Sigma^cM\) defines an operator \(D_{\mathbf L}=c\circ\nabla^{\boldsymbol\Sigma^c}\). If \(F_{\mathbf L}\in\Omega^2(M;i\mathbb R)\) is the curvature of the connection and
\[
 c(F_{\mathbf L})
 :=\frac12\sum_{i,j=1}^{n}F_{\mathbf L}(\mathbf e_i,\mathbf e_j)
 c(\mathbf e^i)c(\mathbf e^j),
\]
then
\begin{equation}
 D_{\mathbf L}^2
 =(\nabla^{\boldsymbol\Sigma^c})^*
 \nabla^{\boldsymbol\Sigma^c}
 +\frac14R_{\mathbf g}I+\frac12c(F_{\mathbf L}).
\label{eq:lichnerowicz-spinc}
\end{equation}
In even dimensions, \(D_{\mathbf L}\) interchanges \(\boldsymbol\Sigma^{c,+}M\) and \(\boldsymbol\Sigma^{c,-}M\).
\end{proposition}

\begin{proof}
The infinitesimal homomorphism of \(\operatorname{Spin}^c(n)\) on spinors is the spin representation on the \(\mathfrak{so}(n)\) component and half the representation of the determinant line on the \(i\mathbb R\) component. Hence,
\[
 \mathbf R^{\boldsymbol\Sigma^c}(\mathbf X,\mathbf Y)
 =\mathbf R^{\boldsymbol\Sigma}(\mathbf X,\mathbf Y)
 +\frac12F_{\mathbf L}(\mathbf X,\mathbf Y)I
\]
in any local spinor trivialization. The general Bochner formula and \eqref{eq:contraccion-curvatura-spinorial} give \eqref{eq:lichnerowicz-spinc}. The graded assertion follows from the fact that the connection preserves chirality and Clifford multiplication is odd.
\end{proof}

\chapter{Sobolev spaces on vector bundles: an intrinsic approach}\label{cap:sobolev-haces}

\begin{semblanzaHistorica}{From functions to sections} Passing from functions to sections does not merely repeat the Euclidean theory component by component. Components depend on the frame, and ordinary derivatives of them produce additional terms under a change of trivialization. A connection organizes these terms and makes it possible to define weak derivatives taking values in a specific tensor bundle. The integration-by-parts formula will bring these two perspectives together: it is computed locally, yet expresses an intrinsic identity and allows the regularity of a section to be recognized without choosing a global frame. \end{semblanzaHistorica}

\section{Weak covariant derivatives and Sobolev spaces on vector bundles}

The preceding construction of Sobolev spaces on manifolds started with smooth sections and completed that space in a covariant norm. There is another description, closer to the Euclidean case: require a section in $L^p(M,\mathbf E)$ to possess weak covariant derivatives through the order under consideration. This chapter develops that formulation on a vector bundle and proves its compatibility with the construction by completion.

For compactly supported smooth sections, the identity \[
\int_M\langle\nabla^{\mathbf{E}} \mathbf{u},\boldsymbol{\Phi}\rangle_{\mathbf{g},\mathbf{h}_{\mathbf{E}}}\,d\lambda_{\mathbf{g}}
=
\int_M\langle \mathbf{u},(\nabla^{\mathbf{E}})^*\boldsymbol{\Phi}\rangle_{\mathbf{h}_{\mathbf{E}}}\,d\lambda_{\mathbf{g}}
\] transfers the derivative from $\mathbf{u}$ to a test section. The right-hand side still makes sense for $\mathbf{u}\in L^p(M,\mathbf{E})$ provided $(\nabla^{\mathbf{E}})^*\boldsymbol{\Phi}$ belongs to the appropriate dual space. This observation leads to the weak covariant derivative and allows elliptic systems on bundles to be formulated without choosing global components.

We work with a smooth vector bundle $\mathbf{E}\longrightarrow M$ equipped with a bundle metric and a smooth connection $\nabla^{\mathbf E}$. On tensor factors constructed from $TM$ and $T^*M$, we use the Levi--Civita connection of $\mathbf g$. The definitions of weak covariant derivative and $W^{m,p}(M,\mathbf E)$ do not require $\nabla^{\mathbf E}$ to be compatible with the bundle metric; whenever compatibility is needed for an integration-by-parts formula, an identification of the formal adjoint, or a subsequent estimate, it will be explicitly included in the hypotheses. The notation $W^{m,p}(M,\mathbf{E})$ denotes the space of sections of $L^p(M,\mathbf{E})$ whose weak covariant derivatives through order $m$ belong to the corresponding $L^p(M,T^{(0,j)}(TM)\otimes\mathbf E)$ spaces of $\mathbf{E}$-valued tensors. This choice retains the notation used in the Euclidean case and treats functions, differential forms, tensor fields, and sections of general bundles uniformly.

To formulate the definition, we first need to construct iterated covariant derivatives and the induced bundle metrics on $T^{(0,s)}(TM)\otimes \mathbf{E}$. We then prove completeness, density of smooth sections, and equivalence with local descriptions. These results give an intrinsic version of the Meyers--Serrin theorem and prepare the embedding, trace, and compactness theorems.

\begin{definition}\label{derivada covariante de orden superior para haces vectoriales}\index{higher-order covariant derivative} Let $M$ be a smooth manifold with or without boundary, let $\mathbf{E}\to M$ be a smooth vector bundle with a connection $\nabla^{\mathbf{E}}$, and let $\nabla^{TM}$ be a connection on $TM$. For $s\in\mathbb{N}$, define the $s$th covariant derivative of $\mathbf{u}\in\Gamma(\mathbf{E})$ as the section $\nabla^s \mathbf{u} \in \Gamma\big(T^{(0,s)}(TM)\otimes \mathbf{E}\big)$ defined recursively as follows: $\nabla^0 \mathbf{u}:=\mathbf{u}$ and, for $s\geq 0$, \[
\nabla^{s+1}\mathbf{u} := \nabla^{T^{(0,s)}(TM)\otimes \mathbf{E}}( \nabla^s \mathbf{u}),
\] where we identify $T^*M\otimes(T^{(0,s)}(TM)\otimes \mathbf{E})$ with $T^{(0,s+1)}(TM)\otimes \mathbf{E}$ by placing the new factor $T^*M$ first. The connection on the right-hand side is the product connection induced by $\nabla^{TM}$ and $\nabla^{\mathbf{E}}$.

Schematically, \[
\begin{tikzcd}[ampersand replacement=\&, column sep=5em, row sep=2em]
 \Gamma(\mathbf{E}) \arrow[r, "\nabla^{\mathbf{E}}"] \&
 \Gamma\bigl(T^{(0,1)}(TM)\otimes \mathbf{E}\bigr)
 \arrow[d, out=0, in=0, looseness=1.4, "{\nabla^{T^{(0,1)}(TM)\otimes \mathbf{E}}}"{right, fill=none}] \\
 %
 \Gamma\bigl(T^{(0,3)}(TM)\otimes \mathbf{E}\bigr)
 \arrow[d, out=180, in=180, looseness=1.4, "{\nabla^{T^{(0,3)}(TM)\otimes \mathbf{E}}}"{left, fill=none}] \&
 \Gamma\bigl(T^{(0,2)}(TM)\otimes \mathbf{E}\bigr) \arrow[l, "{\nabla^{T^{(0,2)}(TM)\otimes \mathbf{E}}}"] \\
 %
 \cdots \arrow[r, "{\nabla^{T^{(0,s-1)}(TM)\otimes \mathbf{E}}}"] \&
 \Gamma\bigl(T^{(0,s)}(TM)\otimes \mathbf{E}\bigr)
\end{tikzcd}
\] \end{definition} Equivalently, \[
(\nabla^{s+1}\mathbf{u})(\mathbf{X}_0,\ldots,\mathbf{X}_s)
:=\bigl(\nabla_{\mathbf{X}_0}(\nabla^s\mathbf{u})\bigr)(\mathbf{X}_1,\ldots,\mathbf{X}_s).
\] In a chart $(x^1,\ldots,x^n)$, with $\boldsymbol{\partial}_i=\partial/\partial x^i$, write \[
(\nabla^s \mathbf{u})_{i_1\cdots i_s}
:=(\nabla^s \mathbf{u})(\boldsymbol{\partial}_{i_1},\ldots,\boldsymbol{\partial}_{i_s}).
\] In particular, \begin{equation}
(\nabla^2\mathbf{u})_{ij}
=\nabla^{\mathbf{E}}_{\boldsymbol{\partial}_i}\bigl(\nabla^{\mathbf{E}}_{\boldsymbol{\partial}_j}\mathbf{u}\bigr)
-\Gamma^k_{ij}\nabla^{\mathbf{E}}_{\boldsymbol{\partial}_k}\mathbf{u}.
\label{eq:componentes-segunda-derivada-convencion-cap9}
\end{equation} This is a component of the tensor $\nabla^2\mathbf{u}$, rather than the mere directional composition: it includes the Christoffel term. At higher orders, we will primarily use $\nabla^s \mathbf{u}$ and its components $(\nabla^s \mathbf{u})_{i_1\cdots i_s}$, according to the same recursion.

The section $\nabla^s\mathbf{u}$ takes values in $T^{(0,s)}(TM)\otimes \mathbf{E}$. To measure it pointwise, we need the bundle metric induced by the Riemannian metric on $M$ and the metric on $\mathbf{E}$.

The same construction defines a bundle metric on every bundle $T^{(k,l)}(TM)\otimes \mathbf{E}$.

\begin{proposition}\label{metrica en tensores con E general} Let $(M,\mathbf{g})$ be a Riemannian manifold (with or without boundary), and let $\mathbf{E}\to M$ be a smooth vector bundle equipped with a bundle metric $\mathbf{h}_{\mathbf{E}}$. Then there is a unique smooth bundle metric $\langle \cdot,\cdot \rangle_{\mathbf{g},\mathbf{h}_{\mathbf{E}}}$ on each $\mathbf{E}$-valued tensor bundle, $T^{(k,l)}(TM)\otimes \mathbf{E}$, such that, for simple tensor fields \[
\boldsymbol{\alpha}_1\otimes\dots\otimes \boldsymbol{\alpha}_k\otimes \boldsymbol{\omega}_1\otimes\dots\otimes \boldsymbol{\omega}_l\otimes \mathbf{u},\quad
\boldsymbol{\beta}_1\otimes\dots\otimes \boldsymbol{\beta}_k\otimes \boldsymbol{\eta}_1\otimes\dots\otimes \boldsymbol{\eta}_l\otimes \mathbf{v},
\] with $\boldsymbol{\alpha}_r,\boldsymbol{\beta}_r\in\mathfrak X(M)$, $\boldsymbol{\omega}_s,\boldsymbol{\eta}_s\in\Omega^1(M)$, and $\mathbf{u},\mathbf{v}\in\Gamma(\mathbf{E})$, we have \[\big\langle \boldsymbol{\alpha}_1\otimes\dots\otimes \boldsymbol{\alpha}_k\otimes \boldsymbol{\omega}_1\otimes\dots\otimes \boldsymbol{\omega}_l\otimes \mathbf{u},
\boldsymbol{\beta}_1\otimes\dots\otimes \boldsymbol{\beta}_k\otimes \boldsymbol{\eta}_1\otimes\dots\otimes \boldsymbol{\eta}_l\otimes \mathbf{v}\big\rangle_{\mathbf{g},\mathbf{h}_{\mathbf{E}}}
=\left(\prod_{s=1}^k\langle \boldsymbol{\alpha}_s,\boldsymbol{\beta}_s\rangle_{\mathbf{g}}\right)\left(\prod_{t=1}^l\langle \boldsymbol{\omega}_t,\boldsymbol{\eta}_t\rangle_{\mathbf{g}}\right)\langle \mathbf{u},\mathbf{v}\rangle_{\mathbf{h}_{\mathbf{E}}}.\]

Moreover, if $(\mathbf{E}_{1},\dots,\mathbf{E}_{n})$ is a smooth local frame of $TM$, with dual coframe $(\boldsymbol{\varepsilon}^{1},\dots,\boldsymbol{\varepsilon}^{n})$, $(\mathbf{e}_1,\dots,\mathbf{e}_{r})$ is a local frame of $\mathbf{E}$, all defined on $U\subseteq M$, and \[
\mathbf{F}=F^{i_1\dots i_k a}_{j_1\dots j_l}\mathbf{E}_{i_1}\otimes\dots\otimes \mathbf{E}_{i_k}\otimes \boldsymbol{\varepsilon}^{j_1}\otimes\dots\otimes \boldsymbol{\varepsilon}^{j_l}\otimes \mathbf{e}_a,
\] \[
\mathbf{G}=G^{r_1\dots r_k b}_{s_1\dots s_l}\mathbf{E}_{r_1}\otimes\dots\otimes \mathbf{E}_{r_k}\otimes \boldsymbol{\varepsilon}^{s_1}\otimes\dots\otimes \boldsymbol{\varepsilon}^{s_l}\otimes \mathbf{e}_b,
\] then \[
\langle \mathbf{F},\mathbf{G}\rangle_{\mathbf{g},\mathbf{h}_{\mathbf{E}}}
=g_{i_1r_1}\cdots g_{i_kr_k}g^{j_1s_1}\cdots g^{j_ls_l}h_{ab}
F^{i_1\dots i_k a}_{j_1\dots j_l}
\overline{G^{r_1\dots r_k b}_{s_1\dots s_l}}.
\] In the real case, the bar is omitted. In the complex case, it appears because the Hermitian convention fixed in Definition~\ref{def:metrica-hermitiana-fibrada} is linear in the first argument and antilinear in the second. When no confusion is possible, we write simply $\langle \mathbf{F},\mathbf{G}\rangle_{\mathbf{g},\mathbf{h}_{\mathbf{E}}}=\langle \mathbf{F},\mathbf{G}\rangle_{\mathbf{h}_{\mathbf{E}}}$. We use the notation $|\mathbf{F}|_{\mathbf{h}_{\mathbf{E}}}:=\sqrt{\langle \mathbf{F},\mathbf{F}\rangle_{\mathbf{h}_{\mathbf{E}}}}$. \end{proposition} \begin{proof} At each point $p\in M$, the metric $\mathbf g_p$ induces the dual inner product $\mathbf g_p^{-1}$ on $T_p^*M$. The tensor product of $k$ copies of $\mathbf g_p$, $l$ copies of $\mathbf g_p^{-1}$, and one copy of $(\mathbf h_{\mathbf E})_p$ defines a positive-definite form on \[
T_p^{(k,l)}(T_pM)\otimes\mathbf E_p.
\] The formula on simple tensors is precisely the definition of this tensor product. Since simple tensors span the fiber, sesquilinearity determines the form uniquely on the entire fiber.

In the frames of the statement, the inner-product matrices of the factors are, respectively, \[
(g_{ir}),\qquad (g^{js}),\qquad (h_{ab}).
\] Expanding $\mathbf F$ and $\mathbf G$ in the tensor basis and using linearity of the Hermitian product in the first variable and antilinearity in the second gives \[
\langle \mathbf{F},\mathbf{G}\rangle_{\mathbf{g},\mathbf{h}_{\mathbf{E}}}
=g_{i_1r_1}\cdots g_{i_kr_k}
g^{j_1s_1}\cdots g^{j_ls_l}h_{ab}
F^{i_1\dots i_k a}_{j_1\dots j_l}
\overline{G^{r_1\dots r_k b}_{s_1\dots s_l}}.
\] In the real case, the same calculation is bilinear and involves no conjugation. The coefficients $g_{ir}$, $g^{js}$, and $h_{ab}$ are smooth in every local frame; hence the resulting family of forms varies smoothly with $p$. Their definition through tensor products shows that expressions in different frames agree on overlaps. This gives a global smooth bundle metric, and the fiberwise uniqueness already proved completes the proof. \end{proof}

With these definitions and conventions, we can define what it means for a section to be an $L^{p}$-section.

\begin{definition}\label{secciones Lp y L1 loc}\index{sections!locally integrable}\index{sections!Lebesgue spaces} Let $(M,\mathbf{g})$ be a Riemannian manifold with or without boundary, let $1\leq p\leq \infty$, and let $\pi\colon \mathbf{E}\longrightarrow M$ be a smooth vector bundle equipped with a bundle metric $\mathbf{h}_{\mathbf{E}}$. We say that a function $\mathbf{u}\colon M\longrightarrow \mathbf{E}$ is an $L^{p}$-section if \begin{enumerate}[label=(\alph*)] \item $\mathbf{u}$ is measurable, that is, for every open subset $U\subseteq \mathbf{E}$ we have $\mathbf{u}^{-1}(U)\in\mathcal{L}_M$. \item $(\pi\circ \mathbf{u})(x)=x$ for almost every $x\in M$. \item $|\mathbf{u}|_{\mathbf{h}_{\mathbf{E}}}\in L^{p}(M,d\lambda_{\mathbf{g}})$. \end{enumerate}

Write \[
\widetilde{L}^{p}(\mathbf{E})
:=
\{\mathbf{u}\colon M\longrightarrow \mathbf{E} \mid \mathbf{u} \text{ is an $L^{p}$ section of } \mathbf{E}\}.
\]

Define an equivalence relation $\sim$ on $\widetilde{L}^{p}(\mathbf{E})$ by \[
\mathbf{u}\sim \mathbf{v}
\quad\Longleftrightarrow\quad
\mathbf{u}(x)=\mathbf{v}(x)\ \text{ for almost every }x\in M.
\]

Finally, define $L^{p}(M,\mathbf{E}):=\widetilde{L}^{p}(\mathbf{E})\big/\sim$. Likewise, define \[
L^{p}_{\mathrm{loc}}(M,\mathbf{E})
:=
\left\{
\mathbf{u}\colon M\longrightarrow \mathbf{E} \ \middle|\
\begin{array}{l}
\mathbf{u} \text{ is measurable},\,
(\pi\circ \mathbf{u})(x)=x \text{ for almost every } x\in M,\\
\forall K\subseteq M \text{ compact},\
\displaystyle\int_{K}|\mathbf{u}|_{\mathbf{h}_{\mathbf{E}}}^{p}\,d\lambda_{\mathbf{g}}<\infty
\end{array}
\right\}\Biggr/\sim,
\] where $\sim$ again denotes equality almost everywhere on $M$. \end{definition}

These spaces can be equipped with norms or seminorms as appropriate: \begin{definition}[$L^{p}$ norms and local seminorms]\label{normas Lp}\index{norms and local seminorms} Let $(M,\mathbf{g})$ be a smooth Riemannian manifold with or without boundary, let $1\le p\le\infty$, and let $\mathbf{E}\to M$ be a smooth vector bundle equipped with a bundle metric $\mathbf{h}_{\mathbf{E}}$.

\begin{enumerate} \item \emph{Global case.} \begin{enumerate} \item If $1\le p<\infty$ and $\mathbf{u}\in L^{p}(M,\mathbf{E})$, define \[
\|\mathbf{u}\|_{L^{p}(M,\mathbf{E})}
:=
\left(\int_M |\mathbf{u}|_{\mathbf{h}_{\mathbf{E}}}^{\,p}\, d\lambda_{\mathbf{g}}\right)^{\frac{1}{p}}.
\] \item If $p=\infty$ and $\mathbf{u}\in L^{\infty}(M,\mathbf{E})$, define \[
\|\mathbf{u}\|_{L^{\infty}(M,\mathbf{E})}
:=
\operatorname*{ess\,sup}_{x\in M}|\mathbf{u}(x)|_{\mathbf{h}_{\mathbf{E}}}.
\] \end{enumerate}

\item \emph{Local case.} \begin{enumerate} \item If $1\le p<\infty$ and $\mathbf{u}\in L^{p}_{\mathrm{loc}}(M,\mathbf{E})$, for each compact set $K\subseteq M$ define the seminorm \[
\|\mathbf{u}\|_{L^{p}(K,\mathbf{E})}
:=
\left(\int_K |\mathbf{u}|_{\mathbf{h}_{\mathbf{E}}}^{\,p}\, d\lambda_{\mathbf{g}}\right)^{\frac{1}{p}}.
\] \item If $p=\infty$ and $\mathbf{u}\in L^{\infty}_{\mathrm{loc}}(M,\mathbf{E})$, for each compact set $K\subseteq M$ define \[
\|\mathbf{u}\|_{L^{\infty}(K,\mathbf{E})}
:=
\operatorname*{ess\,sup}_{x\in K}|\mathbf{u}(x)|_{\mathbf{h}_{\mathbf{E}}}.
\] \end{enumerate} \end{enumerate}

The spaces $L^{p}_{\mathrm{loc}}(M,\mathbf{E})$ are equipped with the family of seminorms $\{\|\cdot\|_{L^{p}(K,\mathbf{E})}\}_{K\subseteq M\ \mathrm{compact}}$. \end{definition}

Our definition of weak derivative uses the formal Hermitian adjoint developed in the preceding chapter. Throughout this chapter, $(\nabla^s)_h^*$ denotes that adjoint; when a previously established formula abbreviates it as $(\nabla^s)^*$, the abbreviation has exactly the same meaning and never denotes the bilinear transpose $(\nabla^s)'$. \begin{definition}\label{def:definicion-sobolev-haces-derivada-variedad-riemanniana-haz-vectorial}\index{higher-order weak covariant derivative@higher-order weak covariant derivative}\index{derivative!weak covariant}\glsadd{derivada-covariante-debil} Let $(M,\mathbf{g})$ be a Riemannian manifold without boundary, and let $\mathbf{E}$ be a smooth vector bundle over $M$ equipped with a bundle metric $\mathbf{h}_{\mathbf{E}}$ and a smooth connection $\nabla^{\mathbf E}$. In constructing $\nabla^s$, use the product connection induced by the Levi--Civita connection of $\mathbf g$ and by $\nabla^{\mathbf E}$. Let $\mathbf{u}\in L_{\text{loc}}^{1}(M,\mathbf{E})$. We say that $\mathbf{v}\in L^{1}_{\text{loc}}(M,T^{(0,s)}(TM)\otimes \mathbf{E})$ is the $s$th weak covariant derivative of $\mathbf{u}$ if, for every $\boldsymbol{\sigma} \in \Gamma_{c}(T^{(0,s)}(TM)\otimes \mathbf{E})$, \[\int_{M}\langle \mathbf{v},\boldsymbol{\sigma} \rangle_{\mathbf{g},\mathbf{h}_{\mathbf{E}}}d\lambda_{\mathbf{g}}=\int_{M}\langle \mathbf{u},(\nabla^{s})_h^{*}\boldsymbol{\sigma} \rangle_{\mathbf{h}_{\mathbf{E}}}d\lambda_{\mathbf{g}}.\] \end{definition}

The weak derivative of a section $\mathbf{u}\in L^{1}_{\text{loc}}(M,\mathbf{E})$ will be unique, as a consequence of a result analogous to Proposition~\ref{prop unicidad derivada debil prop 14.49}. \begin{lemma}[Localization]\label{lema: principio de localizacion para haces vectoriales} Let $(M,\mathbf{g})$ be a Riemannian manifold without boundary, $\mathbf{E}$ a smooth vector bundle with bundle metric $\mathbf{h}_{\mathbf{E}}$, and $\mathbf{v}\in L^{1}_{\mathrm{loc}}(M,\mathbf{E})$ such that \[
\int_{M}\langle \mathbf{v},\boldsymbol{\sigma}\rangle_{\mathbf{h}_{\mathbf{E}}}d\lambda_{\mathbf{g}}=0,\qquad \forall \boldsymbol{\sigma}\in \Gamma_{c}(\mathbf{E}).
\] Then $\mathbf{v}\equiv 0$ almost everywhere on $M$. \end{lemma}

\begin{proof} Consider an open subset $U\subseteq M$ such that $(U,\psi)$ is a smooth chart of $M$ and $U$ is also the domain of a smooth local orthonormal frame of $\mathbf{E}$ over $U$, say $(\mathbf{e}_{1},\dots,\mathbf{e}_{r})$.

Given a section $\boldsymbol{\sigma}\in\Gamma_{c}(\mathbf{E}\restriction_{U})$, extend it by zero outside $U$ to obtain $\widetilde{\boldsymbol{\sigma}}\in\Gamma_{c}(\mathbf{E})$. Applying the hypothesis gives \[
0=\int_M \langle \mathbf{v},\widetilde{\boldsymbol{\sigma}}\rangle_{\mathbf{h}_{\mathbf{E}}}d\lambda_{\mathbf{g}}=\int_U \langle \mathbf{v},\boldsymbol{\sigma}\rangle_{\mathbf{h}_{\mathbf{E}}}d\lambda_{\mathbf{g}}.
\] Since $\mathbf{v}\in L^{1}_{\operatorname{loc}}(M,\mathbf{E})$, we can write $\mathbf{v}=\displaystyle\sum_{i=1}^{r} v^{i} \mathbf{e}_i$ with $v^{i}\in L^{1}_{\mathrm{loc}}(U)$ and $\boldsymbol{\sigma}=\displaystyle\sum_{j=1}^{r} \sigma^{j} \mathbf{e}_j$ with $\sigma^{j}\in C_c^\infty(U)$. Orthonormality of the frame and the convention of linearity in the first variable give \[
0=\int_U \displaystyle\sum_{i=1}^{r} v^{i}\overline{\sigma^{i}}d\lambda_{\mathbf{g}}.
\] Now fix $j\in\{1,\dots,r\}$ and take $\boldsymbol{\sigma}:=f \mathbf{e}_j$ with $f\in C_c^\infty(U)$, obtaining \[
\int_U v^{j}\overline f\,d\lambda_{\mathbf{g}}=0,\qquad \forall f\in C_c^\infty(U).
\]

Passing to local coordinates gives \[
\int_{\psi(U)} (v^{j}\circ\psi^{-1})\overline{(f\circ\psi^{-1})}\bigl(\sqrt{\det(\mathbf{g})}\circ\psi^{-1}\bigr)d\lambda_{n}=0,\qquad \forall f\in C_c^\infty(U).
\] Given arbitrary $\phi\in C_c^\infty(\psi(U))$, define \[
f := \left( \overline{\frac{\phi}{\sqrt{\det(\mathbf{g})}\circ\psi^{-1}}} \right)\circ\psi \in C_c^\infty(U).
\] The function is well defined and smooth because $\sqrt{\det(\mathbf{g})}\circ\psi^{-1}$ is smooth and strictly positive. Moreover, \[
\supp(f)\subseteq
\psi^{-1}\!\left(
\supp\left(\frac{\phi}{\sqrt{\det(\mathbf{g})}\circ\psi^{-1}}\right)
\right),
\], and the set on the right is compact: the support inside $\psi^{-1}$ is compact in $\psi(U)$, and $\psi$ is a homeomorphism. Substituting this choice yields

\[\int_{\psi(U)} (v^{j}\circ\psi^{-1})\phi\,d\lambda_{n}=0,\qquad \forall\phi\in C_c^\infty(\psi(U)).\]

By Proposition~\ref{prop unicidad derivada debil prop 14.49}, we conclude that $v^{j}\circ\psi^{-1}=0$ almost everywhere on $\psi(U)$, and hence $v^{j}=0$ almost everywhere on $U$. Since $j\in\{1,\dots,r\}$ was arbitrary, it follows that $\mathbf{v}=0$ almost everywhere on $U$.

Finally, since $M$ is Lindelöf, it can be covered by countably many open subsets of this kind, on each of which $\mathbf{v}=0$ almost everywhere. It follows that $\mathbf{v}\equiv 0$ almost everywhere on $M$. \end{proof}

\begin{corollary}\label{cor:definicion-sobolev-haces-esima-derivada-covariante-debil-esta-unica} Let $(M,\mathbf{g})$ be a Riemannian manifold without boundary, and let $\mathbf{E}\to M$ be a smooth vector bundle with a bundle metric. If $\mathbf{u}\in L_{\operatorname{loc}}^{1}(M,\mathbf{E})$ has an $s$th weak covariant derivative, it is unique, and we denote it by $\nabla ^{s}_{w}\mathbf{u}$, where $w$ stands for ``weak.'' \end{corollary} \begin{proof} Suppose $\mathbf{v}_1$ and $\mathbf{v}_2$ are two $s$th weak covariant derivatives of $\mathbf{u}$. For every $\boldsymbol{\sigma}\in\Gamma_c(T^{(0,s)}(TM)\otimes\mathbf{E})$, their two defining identities have the same right-hand side. Subtraction gives \[
\int_M
\left\langle\mathbf{v}_1-\mathbf{v}_2,\boldsymbol{\sigma}\right\rangle_{\mathbf{g},\mathbf{h}_{\mathbf{E}}}
\,d\lambda_{\mathbf{g}}=0.
\] Lemma~\ref{lema: principio de localizacion para haces vectoriales}, applied to the bundle $T^{(0,s)}(TM)\otimes\mathbf{E}$ with its induced bundle metric, implies that $\mathbf{v}_1-\mathbf{v}_2=0$ almost everywhere. Therefore $\mathbf{v}_1=\mathbf{v}_2$ almost everywhere. \end{proof}

The weak covariant derivative is a local operator, just like its classical counterpart: \begin{proposition}\label{prop: derivada covariante debil es operador local} Let $(M,\mathbf{g})$ be a Riemannian manifold without boundary and $\mathbf{E}\longrightarrow M$ a smooth vector bundle with bundle metric $\mathbf{h}_{\mathbf{E}}$. Let $\mathbf{u}\in L^{1}_{\operatorname{loc}}(M,\mathbf{E})$ be such that $\mathbf{u}$ has an $s$th weak covariant derivative $\nabla^{s}_{w}\mathbf{u}\in L^{1}_{\operatorname{loc}}(M,T^{(0,s)}(TM)\otimes \mathbf{E})$. Let $U\subseteq M$ be open. If $\mathbf{u}\equiv 0$ almost everywhere on $U$, then $\nabla^{s}_{w}\mathbf{u}\equiv 0$ almost everywhere on $U$. \end{proposition} \begin{proof} Let $\boldsymbol{\sigma} \in \Gamma_{c}(T^{(0,s)}(TU)\otimes \mathbf{E}\restriction_{U})$. We can define $\widetilde{\boldsymbol{\sigma} }\in \Gamma_{c}(T^{(0,s)}(TM)\otimes \mathbf{E})$ by \[\widetilde{\boldsymbol{\sigma}}(p)=\begin{cases}
 \boldsymbol{\sigma}(p)&\text{if $p\in U$}\\
 0_{\mathbf{E}_{p}}&\text{if $p\notin U$}
 \end{cases}.\] Since $U$ is open and $\supp(\boldsymbol{\sigma} )\subseteq U$, $\widetilde{\boldsymbol{\sigma}}$ is smooth, so $\widetilde{\boldsymbol{\sigma}}\in \Gamma_{c}(T^{(0,s)}(TM)\otimes \mathbf{E})$. By the definition of weak covariant derivative, we have \[\int_{M}\langle \nabla_{w}^{s}\mathbf{u},\widetilde{\boldsymbol{\sigma}}\rangle _{\mathbf{g},\mathbf{h}_{\mathbf{E}}}d\lambda_{\mathbf{g}}=\int_{M}\langle \mathbf{u},(\nabla^{s})^{*}\widetilde{\boldsymbol{\sigma}}\rangle_{\mathbf{h}_{\mathbf{E}}}d\lambda_{\mathbf{g}}.\] Since $\supp(\widetilde{\boldsymbol{\sigma}})=\supp(\boldsymbol{\sigma})\subseteq U$ and by Proposition~\ref{pdo es local} ($(\nabla^{s})^{*}$ is a local operator), we have $\supp((\nabla^{s})^{*}\widetilde{\boldsymbol{\sigma}})\subseteq \supp(\widetilde{\boldsymbol{\sigma}})\subseteq U$. Then \[\int_{U}\langle \nabla_{w}^{s}\mathbf{u},\widetilde{\boldsymbol{\sigma}}\rangle _{\mathbf{g},\mathbf{h}_{\mathbf{E}}}d\lambda_{\mathbf{g}}=\int_{M}\langle \nabla_{w}^{s}\mathbf{u},\widetilde{\boldsymbol{\sigma}}\rangle _{\mathbf{g},\mathbf{h}_{\mathbf{E}}}d\lambda_{\mathbf{g}}=\int_{M}\langle \mathbf{u},(\nabla^{s})^{*}\widetilde{\boldsymbol{\sigma}}\rangle_{\mathbf{h}_{\mathbf{E}}}d\lambda_{\mathbf{g}}=\int_{U}\langle \mathbf{u},(\nabla^{s})^{*}\widetilde{\boldsymbol{\sigma}}\rangle_{\mathbf{h}_{\mathbf{E}}}d\lambda_{\mathbf{g}}=0,\] where the last equality holds because $\mathbf{u}\equiv 0$ on $U$. Therefore Lemma~\ref{lema: principio de localizacion para haces vectoriales} gives $\nabla_{w}^{s}\mathbf{u}\equiv 0$ on $U$. \end{proof}

With these definitions, we can define Sobolev spaces by weak covariant derivatives on any smooth vector bundle. \begin{definition}[Global and local Sobolev spaces]\label{Sobolev global y local}\index{global and local Sobolev spaces@global and local Sobolev spaces} Let $(M,\mathbf{g})$ be a Riemannian manifold without boundary, and let $\pi\colon \mathbf{E}\longrightarrow M$ be a smooth vector bundle equipped with a bundle metric $\mathbf{h}_{\mathbf{E}}$ and a smooth connection $\nabla^{\mathbf E}$. Let $m\in\mathbb{N}_0$ and $1\le p\le\infty$.

\begin{enumerate} \item \emph{Global Sobolev spaces.} Define \[
W^{m,p}(M,\mathbf{E})
:=
\left\{
\mathbf{u}\in L^{1}_{\mathrm{loc}}(M,\mathbf{E})
\;\middle|\;
\nabla_w^{s}\mathbf{u}\in L^{p}\big(M,T^{(0,s)}(TM)\otimes \mathbf{E}\big)
\ \text{for every } s\in\{0,\dots,m\}
\right\}.
\] The space $W^{m,p}(M,\mathbf{E})$ is a Banach space with norm \[
\|\mathbf{u}\|_{W^{m,p}(M,\mathbf{E})}
:=
\begin{cases}
\displaystyle
\left(
\displaystyle\sum_{s=0}^{m}
\int_M |\nabla_w^{s}\mathbf{u}|_{\mathbf{g},\mathbf{h}_{\mathbf{E}}}^{\,p}\, d\lambda_{\mathbf{g}}
\right)^{\frac{1}{p}},
& 1\le p<\infty, \\[1.2em]
\displaystyle
\max_{0\le s\le m}
\operatorname*{ess\,sup}_{x\in M}
|\nabla_w^{s}\mathbf{u}(x)|_{\mathbf{g},\mathbf{h}_{\mathbf{E}}},
& p=\infty.
\end{cases}
\]

\item \emph{Local Sobolev spaces.} Define \[
W^{m,p}_{\mathrm{loc}}(M,\mathbf{E})
:=
\left\{
\mathbf{u}\in L^{1}_{\mathrm{loc}}(M,\mathbf{E})
\;\middle|\;
\nabla_w^{s}\mathbf{u}\in L^{p}_{\mathrm{loc}}\big(M,T^{(0,s)}(TM)\otimes \mathbf{E}\big)
\ \text{for every } s\in\{0,\dots,m\}
\right\}.
\] For each compact set $K\subseteq M$, equip $W^{m,p}_{\mathrm{loc}}(M,\mathbf{E})$ with the local \emph{seminorm} \[
\|\mathbf{u}\|_{W^{m,p}(M,\mathbf{E};K)}
:=
\begin{cases}
\displaystyle
\left(
\displaystyle\sum_{s=0}^{m}
\int_K |\nabla_w^{s}\mathbf{u}|_{\mathbf{g},\mathbf{h}_{\mathbf{E}}}^{\,p}\, d\lambda_{\mathbf{g}}
\right)^{\frac{1}{p}},
& 1\le p<\infty, \\[1.2em]
\displaystyle
\max_{0\le s\le m}
\operatorname*{ess\,sup}_{x\in K}
|\nabla_w^{s}\mathbf{u}(x)|_{\mathbf{g},\mathbf{h}_{\mathbf{E}}},
& p=\infty.
\end{cases}
\] The space $W^{m,p}_{\mathrm{loc}}(M,\mathbf{E})$ is equipped with the family of seminorms $\{\|\cdot\|_{W^{m,p}(M,\mathbf{E};K)}\}_{K\subseteq M\ \mathrm{compact}}$. \end{enumerate} \end{definition}

\label{conv:notacion-sobolev-conexion}
When comparing different connections on the same bundle, we write
$W^{m,p}(M,\mathbf E;\nabla^{\mathbf E})$ and
$\|\mathbf u\|_{W^{m,p}(M,\mathbf E;\nabla^{\mathbf E})}$ to specify
the datum used in Definition~\ref{Sobolev global y local}.
If the metrics also vary, these data are written after the
semicolon: $W^{m,p}(M,\mathbf E;\mathbf g,\mathbf h_{\mathbf E},
\nabla^{\mathbf E})$. For a tensor bundle over $M$, the notation
$W^{m,p}(M,\mathbf E;\mathbf g)$ includes the tensor metric and the
connection induced by $\mathbf g$. When the data are fixed, we retain
the notation $W^{m,p}(M,\mathbf E)$.

\begin{exercise}\label{ejer:definicion-sobolev-haces-pruebe-espacio-normado} Under the hypotheses of the preceding definition, where $M$ is a smooth manifold with or without boundary, prove that $(W^{m,p}(M,\mathbf{E}),\|\cdot\|_{W^{m,p}(M,\mathbf{E})})$ is a normed space. \end{exercise} \begin{remark}\label{obs:definicion-sobolev-haces-caso-trivial-isomorfismo-natural-haces-fibra} In the trivial case $\mathbf{E}=M\times \mathbb{R}$, there is a natural bundle isomorphism \[
T^{(0,s)}(TM)\otimes (M\times \mathbb{R}) \cong T^{(0,s)}(TM),
\] given fiberwise by $v\otimes \alpha \mapsto \alpha v$. (More explicitly, considering $v_p\in T^{(0,s)}(T_pM)$, $(p,\alpha) \in \{p\}\times \mathbb{R}, $, the assignment is $v_p\otimes (p,\alpha)\mapsto \alpha v_p\in T^{(0,s)}(T_pM)$.) Under this identification, the Sobolev spaces reduce to the usual spaces for functions \[
W^{m,p}(M):=\Big\{u\in L^{1}_{\mathrm{loc}}(M) \mid \nabla^{s}_{w}u\in L^{p}\big(M,T^{(0,s)}(TM)\big)\text{for every }s\in\{0,\dots,m\}\Big\}.
\] \end{remark} \subsection{Fundamental properties: uniform convexity, reflexivity, and separability} We now prove the fundamental properties of Sobolev spaces on vector bundles, which will be analogous to their Euclidean counterparts. First, we show that $W^{m,p}(M,\mathbf{E})$ is a Banach space. The following lemma is particularly useful for this purpose: \begin{lemma}\label{derivada debil es operador cerrado} Let $(M,\mathbf{g})$ be a Riemannian manifold without boundary, and let $\mathbf{E}\longrightarrow M$ be a smooth vector bundle with a bundle metric $\mathbf{h}_{\mathbf{E}}$ and a connection $\nabla^{\mathbf{E}}$. Let $1\leq p\leq \infty$. If $(\mathbf{u}_{k})_{k\in \mathbb{N}}\subseteq W^{m,p}(M,\mathbf{E})$ is such that $\mathbf{u}_{k}\to \mathbf{u}$ in $L^{p}(M,\mathbf{E})$ and $\nabla^{j}_{w}\mathbf{u}_{k}\to \mathbf{v}_{j}$ in $L^{p}(M,T^{(0,j)}(TM)\otimes \mathbf{E})$ as $k\to \infty$ for every $j\in \{1,\dots,m\}$, then $\mathbf{u}\in W^{m,p}(M,\mathbf{E})$ and $\nabla_{w}^{j}\mathbf{u}=\mathbf{v}_{j}$ for every $j\in\{1,\dots,m\}$. \end{lemma} \begin{proof} Let $j\in\{1,\dots,m\}$ be arbitrary and let $\boldsymbol{\sigma}\in\Gamma_c\big(T^{(0,j)}(TM)\otimes \mathbf{E}\big)$. Since $\mathbf{u}_k\in W^{m,p}(M,\mathbf{E})$, the definition of weak covariant derivative of order $j$ gives \[
\int_M \langle \nabla_w^{j}\mathbf{u}_k,\boldsymbol{\sigma}\rangle_{\mathbf{g},\mathbf{h}_{\mathbf{E}}}\, d\lambda_{\mathbf{g}}
=
\int_M \langle \mathbf{u}_k,(\nabla^{j})^{*}\boldsymbol{\sigma}\rangle_{\mathbf{h}_{\mathbf{E}}}\, d\lambda_{\mathbf{g}}.
\] By Hölder's inequality (Proposition~\ref{desigualdad de holder}), with $p'$ conjugate to $p$, we have \[
\left|
\int_M \langle \mathbf{u}_k-\mathbf{u},(\nabla^{j})^{*}\boldsymbol{\sigma}\rangle_{\mathbf{h}_{\mathbf{E}}}\, d\lambda_{\mathbf{g}}
\right|
\le
\|\mathbf{u}_k-\mathbf{u}\|_{L^p(M,\mathbf{E})}\,
\|(\nabla^{j})^{*}\boldsymbol{\sigma}\|_{L^{p'}(M,\mathbf{E})},
\] and, since $\mathbf{u}_k\to \mathbf{u}$ in $L^p(M,\mathbf{E})$, it follows that \[
\int_M \langle \mathbf{u}_k,(\nabla^{j})^{*}\boldsymbol{\sigma}\rangle_{\mathbf{h}_{\mathbf{E}}}\, d\lambda_{\mathbf{g}}
\longrightarrow
\int_M \langle \mathbf{u},(\nabla^{j})^{*}\boldsymbol{\sigma}\rangle_{\mathbf{h}_{\mathbf{E}}}\, d\lambda_{\mathbf{g}}.
\] Likewise, Hölder's inequality (Proposition~\ref{desigualdad de holder}) gives \[
\left|
\int_M \langle \nabla_w^{j}\mathbf{u}_k-\mathbf{v}_j,\boldsymbol{\sigma}\rangle_{\mathbf{g},\mathbf{h}_{\mathbf{E}}}\, d\lambda_{\mathbf{g}}
\right|
\le
\|\nabla_w^{j}\mathbf{u}_k-\mathbf{v}_j\|_{L^p(M,T^{(0,j)}(TM)\otimes \mathbf{E})}\,
\|\boldsymbol{\sigma}\|_{L^{p'}(M,T^{(0,j)}(TM)\otimes \mathbf{E})},
\] and, since $\nabla_w^{j}\mathbf{u}_k\to \mathbf{v}_j$ in $L^p(M,T^{(0,j)}(TM)\otimes \mathbf{E})$, we obtain \[
\int_M \langle \nabla_w^{j}\mathbf{u}_k,\boldsymbol{\sigma}\rangle_{\mathbf{g},\mathbf{h}_{\mathbf{E}}}\, d\lambda_{\mathbf{g}}
\longrightarrow
\int_M \langle \mathbf{v}_j,\boldsymbol{\sigma}\rangle_{\mathbf{g},\mathbf{h}_{\mathbf{E}}}\, d\lambda_{\mathbf{g}}.
\] Passing to the limit in the initial identity yields \[
\int_M \langle \mathbf{v}_j,\boldsymbol{\sigma}\rangle_{\mathbf{g},\mathbf{h}_{\mathbf{E}}}\, d\lambda_{\mathbf{g}}
=
\int_M \langle \mathbf{u},(\nabla^{j})^{*}\boldsymbol{\sigma}\rangle_{\mathbf{h}_{\mathbf{E}}}\, d\lambda_{\mathbf{g}}
\qquad
\forall\,\boldsymbol{\sigma}\in\Gamma_c\big(T^{(0,j)}(TM)\otimes \mathbf{E}\big).
\] By the definition of weak covariant derivative of order $j$, this implies that $\mathbf{u}$ admits a weak covariant derivative of order $j$ and that $\nabla_w^{j}\mathbf{u}=\mathbf{v}_j$.

Since $j\in\{1,\dots,m\}$ was arbitrary, we conclude that $\mathbf{u}\in W^{m,p}(M,\mathbf{E})$ and $\nabla_w^{j}\mathbf{u}=\mathbf{v}_j$ for every $j\in\{1,\dots,m\}$. \end{proof}

To show that the spaces $W^{m,p}(M,\mathbf{E})$ are Banach, we will, as the reader may expect, use the preceding lemma; first, however, we need to show that $L^{p}(M,\mathbf{E})$ is a Banach space. \begin{theorem}\label{teo:Banach-LpE-total}\index{Lebesgue space@Lebesgue space!sections of a vector bundle as a Banach space} Let $(M,\mathbf{g})$ be a Riemannian manifold without boundary, and let $\pi_{\mathbf{E}}\colon \mathbf{E}\longrightarrow M$ be a smooth vector bundle equipped with a bundle metric $\mathbf{h}_{\mathbf{E}}$. If $1\le p\le \infty$, then $L^{p}(M,\mathbf{E})$ is a Banach space. \end{theorem}

\begin{proof} Let $(\mathbf{s}_n)_{n\in\mathbb N}$ be a Cauchy sequence in $L^{p}(M,\mathbf{E})$.

First consider the case $1\le p<\infty$. There is a subsequence $(\mathbf{s}_{n_k})_{k\in\mathbb N}$ such that \[
\|\mathbf{s}_{n_{k+1}}-\mathbf{s}_{n_k}\|_{L^{p}(M,\mathbf{E})}<\frac{1}{2^k}
\qquad \text{for every }k\in\mathbb N.
\]

Define \[
G_K(x):=\displaystyle\sum_{k=1}^{K}|\mathbf{s}_{n_{k+1}}(x)-\mathbf{s}_{n_k}(x)|_{\mathbf{h}_{\mathbf{E}}}.
\] Then $(G_K)$ is increasing and, by Minkowski's inequality in $L^p(M)$ (Theorem~\ref{teo:b6-espacios-lp-espacios-de-lebesgue}), \[
\|G_K\|_{L^{p}(M)}
\le
\displaystyle\sum_{k=1}^{K}
\|\mathbf{s}_{n_{k+1}}-\mathbf{s}_{n_k}\|_{L^{p}(M,\mathbf{E})}
<
1.
\] By Theorem~\ref{convergencia monotona}, there exists $G\in L^{p}(M)$ such that $G_K\to G$ pointwise and $G(x)<\infty$ almost everywhere.

For each $k$, there exists a null set $Z_k\subset M$ such that $\pi_{\mathbf{E}}(\mathbf{s}_{n_k}(x))=x$ for every $x\in M\setminus Z_k$. Define \[
Z:=\{x\in M\mid G(x)=\infty\}\cup\bigcup_{k=1}^{\infty}Z_k.
\] Then $Z$ is a null set, and for every $x\in M\setminus Z$ we have $G(x)<\infty$ and $\pi_{\mathbf{E}}(\mathbf{s}_{n_k}(x))=x$ for every $k$.

For $x\in M\setminus Z$, the series \[
\displaystyle\sum_{k=1}^{\infty}(\mathbf{s}_{n_{k+1}}(x)-\mathbf{s}_{n_k}(x))
\] converges absolutely in the fiber $\mathbf{E}_x$, since the sum of the norms is $G(x)<\infty$. Since $\mathbf{E}_x$ is finite dimensional, it is complete, and the series converges to a well-defined vector in $\mathbf{E}_x$.

Define \[
\mathbf{s}(x)
:=
\begin{cases}
\displaystyle
\mathbf{s}_{n_1}(x)
+
\displaystyle\sum_{k=1}^{\infty}
(\mathbf{s}_{n_{k+1}}(x)-\mathbf{s}_{n_k}(x)),
& x\in M\setminus Z, \\[6pt]
0_x, & x\in Z,
\end{cases}
\] where $0_x$ denotes the zero vector of the fiber $\mathbf{E}_x$.

For $x\in M\setminus Z$, the sum telescopes, giving \[
\mathbf{s}_{n_{k+1}}(x)\to \mathbf{s}(x)
\qquad \text{a.e. on }M.
\] Each $\mathbf{s}_{n_k}$ is measurable, so $\mathbf{s}$ is measurable as the almost-everywhere pointwise limit of measurable functions. Moreover, by construction $\pi_{\mathbf{E}}(\mathbf{s}(x))=x$ for every $x\in M$, so $\mathbf{s}$ is a measurable section of $\mathbf{E}$.

Let $\varepsilon>0$. Since $(\mathbf{s}_n)$ is Cauchy, there exists $N$ such that \[
\|\mathbf{s}_n-\mathbf{s}_m\|_{L^{p}(M,\mathbf{E})}<\varepsilon
\qquad \text{for every }n,m\ge N.
\] For fixed $n\ge N$, we have \[
|\mathbf{s}_n(x)-\mathbf{s}_{n_k}(x)|_{\mathbf{h}_{\mathbf{E}}}^{p}
\to
|\mathbf{s}_n(x)-\mathbf{s}(x)|_{\mathbf{h}_{\mathbf{E}}}^{p}
\qquad \text{a.e.}
\] Since $n_k\longrightarrow\infty$, there exists $k_0\in\mathbb N$ such that $n_k\geq N$ for every $k\geq k_0$. Applying Fatou's lemma (Lemma~\ref{lem:fatou}), since the functions are nonnegative and measurable, \[
\int_M |\mathbf{s}_n-\mathbf{s}|_{\mathbf{h}_{\mathbf{E}}}^{p}\,d\lambda_{\mathbf{g}}
\le
\liminf_{k\to\infty}
\int_M |\mathbf{s}_n-\mathbf{s}_{n_k}|_{\mathbf{h}_{\mathbf{E}}}^{p}\,d\lambda_{\mathbf{g}}
\le
\varepsilon^{p}.
\] Therefore $\|\mathbf{s}_n-\mathbf{s}\|_{L^{p}(M,\mathbf{E})}\le\varepsilon$ for every $n\ge N$, and we conclude that $\mathbf{s}_n\to \mathbf{s}$ in $L^{p}(M,\mathbf{E})$. Taking a fixed index $n_0$, \[
\|\mathbf{s}\|_{L^{p}(M,\mathbf{E})}
\le
\|\mathbf{s}_{n_0}\|_{L^{p}(M,\mathbf{E})}+\|\mathbf{s}-\mathbf{s}_{n_0}\|_{L^{p}(M,\mathbf{E})}<\infty,
\] so $\mathbf{s}\in L^{p}(M,\mathbf{E})$.

Now consider the case $p=\infty$. By the definition of the essential supremum norm, for each pair $n,m\in\mathbb N$ there is a null set $Z_{n,m}\subset M$ such that \[
|\mathbf{s}_n(x)-\mathbf{s}_m(x)|_{\mathbf{h}_{\mathbf{E}}}\le \|\mathbf{s}_n-\mathbf{s}_m\|_{L^\infty(M,\mathbf{E})}
\qquad \text{for every }x\in M\setminus Z_{n,m}.
\] In addition, for each $n\in\mathbb N$ there is a null set $W_n\subset M$ such that $\pi_{\mathbf{E}}(\mathbf{s}_n(x))=x$ for every $x\in M\setminus W_n$. Define \[
Z:=\left(\bigcup_{n,m=1}^{\infty} Z_{n,m}\right)\cup\left(\bigcup_{n=1}^{\infty}W_n\right).
\] Then $Z$ is a null set and, for every $x\in M\setminus Z$ and every $n,m\in\mathbb N$, \[
|\mathbf{s}_n(x)-\mathbf{s}_m(x)|_{\mathbf{h}_{\mathbf{E}}}\le \|\mathbf{s}_n-\mathbf{s}_m\|_{L^\infty(M,\mathbf{E})}.
\]

Since $(\mathbf{s}_n)$ is Cauchy in $L^\infty(M,\mathbf{E})$, the preceding inequality implies that, for each $x\in M\setminus Z$, the sequence $(\mathbf{s}_n(x))_{n\in\mathbb N}$ is Cauchy in $\mathbf{E}_x$. Since $\mathbf{E}_x$ is finite dimensional, it is complete, and there exists $\mathbf{s}(x)\in \mathbf{E}_x$ such that \[
\mathbf{s}_n(x)\longrightarrow \mathbf{s}(x)
\qquad \text{as }n\to\infty.
\] Define \[
\mathbf{s}(x):=
\begin{cases}
\displaystyle \lim_{n\to\infty}\mathbf{s}_n(x), & x\in M\setminus Z, \\[6pt]
0_x, & x\in Z,
\end{cases}
\] where $0_x$ denotes the zero vector of the fiber $\mathbf{E}_x$. By construction, $\pi_{\mathbf{E}}(\mathbf{s}(x))=x$ for every $x\in M\setminus Z$, and also for $x\in Z$ because $\pi_{\mathbf{E}}(0_x)=x$. Moreover, $\mathbf{s}$ is measurable as the almost-everywhere pointwise limit of the measurable sections $\mathbf{s}_n$.

Let $\varepsilon>0$. Since $(\mathbf{s}_n)$ is Cauchy in $L^\infty(M,\mathbf{E})$, there exists $N\in\mathbb N$ such that \[
\|\mathbf{s}_n-\mathbf{s}_m\|_{L^\infty(M,\mathbf{E})}<\varepsilon
\qquad \text{for every }n,m\ge N.
\] By the construction of $Z$, for every $x\in M\setminus Z$ and every $n,m\ge N$ we have \[
|\mathbf{s}_n(x)-\mathbf{s}_m(x)|_{\mathbf{h}_{\mathbf{E}}}\le \|\mathbf{s}_n-\mathbf{s}_m\|_{L^\infty(M,\mathbf{E})}<\varepsilon.
\] For fixed $n\ge N$, taking the limit as $m\to\infty$ and using $\mathbf{s}_m(x)\to \mathbf{s}(x)$ for $x\in M\setminus Z$ gives \[
|\mathbf{s}_n(x)-\mathbf{s}(x)|_{\mathbf{h}_{\mathbf{E}}}\le \varepsilon
\qquad \text{for every }x\in M\setminus Z.
\] Therefore \[
\|\mathbf{s}_n-\mathbf{s}\|_{L^\infty(M,\mathbf{E})}\le \varepsilon
\qquad \text{for every }n\ge N,
\], and we conclude that $\mathbf{s}_n\to \mathbf{s}$ in $L^\infty(M,\mathbf{E})$. In particular, taking $n=N$ gives $\|\mathbf{s}-\mathbf{s}_N\|_{L^\infty(M,\mathbf{E})}\le\varepsilon$, and since $\mathbf{s}_N\in L^\infty(M,\mathbf{E})$, we deduce that $\mathbf{s}\in L^\infty(M,\mathbf{E})$.

\medskip We conclude that every Cauchy sequence in $L^{p}(M,\mathbf{E})$ converges to an element of $L^{p}(M,\mathbf{E})$. Therefore $L^{p}(M,\mathbf{E})$ is complete for every $1\le p\le\infty$. \end{proof}

A section together with its weak covariant derivatives can be viewed as an element of the product of spaces of integrable sections specified below. The resulting embedding is isometric and transfers functional properties of the product to the Sobolev space. \begin{theorem}[Isometric embedding by covariant derivatives]\label{encaje isometrico en Lp sobolev haces}\index{isometric embedding@isometric embedding} Let $(M,\mathbf{g})$ be a Riemannian manifold without boundary, and let $\mathbf{E}\longrightarrow M$ be a smooth vector bundle with a bundle metric $\mathbf{h}_{\mathbf{E}}$ and a smooth connection $\nabla^{\mathbf{E}}$. Let $m\in\mathbb N_0$ and $1\leq p\leq \infty$. Then the operator $J_{m}\colon W^{m,p}(M,\mathbf{E})\subseteq L^{p}(M,\mathbf{E})\longrightarrow \displaystyle\prod_{j=0}^{m}L^{p}(M,T^{(0,j)}(TM)\otimes \mathbf{E})$ given by $J_{m}(\mathbf{u})=\left(\nabla_w^{j}\mathbf{u}\right)_{j=0}^{m}$ is an isometry onto its image. Moreover, $J_{m}(W^{m,p}(M,\mathbf{E}))$ is closed, and therefore $W^{m,p}(M,\mathbf{E})$ is a Banach space. \end{theorem} \begin{proof} By the definition of the Sobolev norm, if $1\le p<\infty$ we have \[
\|\mathbf{u}\|_{W^{m,p}(M,\mathbf{E})}=
\left(\displaystyle\sum_{j=0}^{m}\|\nabla_w^{j}\mathbf{u}\|_{L^p(M,T^{(0,j)}(TM)\otimes \mathbf{E})}^{p}\right)^{\frac{1}{p}},
\] whereas for $p=\infty$ we have \[
\|\mathbf{u}\|_{W^{m,\infty}(M,\mathbf{E})}=\max_{0\le j\le m}\|\nabla_w^{j}\mathbf{u}\|_{L^\infty(M,T^{(0,j)}(TM)\otimes \mathbf{E})}.
\] Equip the product $\displaystyle\prod_{j=0}^{m}L^{p}(M,T^{(0,j)}(TM)\otimes \mathbf{E})$ with the $\ell^{p}(\{0,\ldots,m\})$ norm if $1\le p<\infty$ and the maximum norm if $p=\infty$. Then, for every $\mathbf{u}\in W^{m,p}(M,\mathbf{E})$, \[
\|J_m(\mathbf{u})\|_{\ell^p(\{0,\ldots,m\})}=
\left\|(\nabla_w^{j}\mathbf{u})_{j=0}^{m}\right\|_{\ell^p(\{0,\ldots,m\})}
=
\|\mathbf{u}\|_{W^{m,p}(M,\mathbf{E})},
\] showing that $J_m$ is an isometry onto its image.

Now prove that $J_m(W^{m,p}(M,\mathbf{E}))$ is closed. Let \[
J_m(\mathbf{u}_k)=\big(\nabla_w^{j}\mathbf{u}_k\big)_{j=0}^{m}\in J_m(W^{m,p}(M,\mathbf{E}))
\] be a sequence such that \[
J_m(\mathbf{u}_k)\longrightarrow (\mathbf{w}_j)_{j=0}^{m}
\quad\text{in}\quad
\prod_{j=0}^{m}L^{p}(M,T^{(0,j)}(TM)\otimes \mathbf{E}).
\] In particular, \[
\mathbf{u}_k\to \mathbf{w}_0 \ \text{in } L^p(M,\mathbf{E}),
\qquad
\nabla_w^{j}\mathbf{u}_k\to \mathbf{w}_j \ \text{in } L^p\big(M,T^{(0,j)}(TM)\otimes \mathbf{E}\big),
\quad j=1,\dots,m.
\] Applying Lemma~\ref{derivada debil es operador cerrado} gives $\mathbf{w}_0\in W^{m,p}(M,\mathbf{E})$ and \[
\nabla_w^{j}\mathbf{w}_0=\mathbf{w}_j,
\qquad j=1,\dots,m.
\] Therefore \[
(\mathbf{w}_j)_{j=0}^{m}=J_m(\mathbf{w}_0)\in J_m(W^{m,p}(M,\mathbf{E})),
\] proving that $J_m(W^{m,p}(M,\mathbf{E}))$ is closed.

Finally, since $\displaystyle\prod_{j=0}^{m}L^{p}(M,T^{(0,j)}(TM)\otimes \mathbf{E})$ is a Banach space and $J_m(W^{m,p}(M,\mathbf{E}))$ is a closed subspace, $J_m(W^{m,p}(M,\mathbf{E}))$ is a Banach space. Since $J_m\colon W^{m,p}(M,\mathbf{E})\longrightarrow J_m(W^{m,p}(M,\mathbf{E}))$ is an isometry onto its image, we conclude that $W^{m,p}(M,\mathbf{E})$ is a Banach space. \end{proof}

Among the properties of Sobolev spaces are, in certain cases, uniform convexity, reflexivity, and separability. We record uniform convexity in the following definition:

\begin{definition}[Uniform convexity]\label{def convexidad uniforme}\index{uniform convexity} Let $(X,\|\cdot\|)$ be a normed space. We say that $X$ is \emph{uniformly convex} if for every $\varepsilon>0$ there exists $\delta>0$ such that, for any $x,y\in X$ with \[
\|x\|=\|y\|=1
\quad\text{and}\quad
\|x-y\|\ge\varepsilon,
\], we have \[
\left\|\frac{x+y}{2}\right\|\le 1-\delta.
\] \end{definition} To establish uniform convexity, we first recall an important theorem of functional analysis, the Milman--Pettis theorem: \begin{theorem}[Milman--Pettis]\label{teorema milman pettis}\index{Milman Pettis@Milman--Pettis} Every uniformly convex Banach space is reflexive. \end{theorem} The following inequalities are named after James A. Clarkson \cite{Clarkson1936}, who studied bounds of this kind for norms on Lebesgue spaces. In the special case of a Hilbert space, its geometric structure considerably simplifies the argument, allowing these inequalities to be deduced directly from the parallelogram law.

\begin{lemma}[Clarkson-type inequalities in Hilbert spaces]\label{lem:clarkson-hilbert}\index{Clarkson-type inequalities in Hilbert spaces@Clarkson-type inequalities in Hilbert spaces} Let $(H,\langle\cdot,\cdot\rangle)$ be a Hilbert space and $\|\cdot\|$ its induced norm. \begin{enumerate} \item If $2\le p<\infty$, then for any $u,v\in H$, \[
\left\|\frac{u+v}{2}\right\|^{p}
+
\left\|\frac{u-v}{2}\right\|^{p}
\le
\frac{\|u\|^{p}+\|v\|^{p}}{2}.
\]

\item If $1<p\le 2$ and $q$ is the conjugate exponent of $p$, then for any $u,v\in H$, \[
\left\|\frac{u+v}{2}\right\|^{q}
+
\left\|\frac{u-v}{2}\right\|^{q}
\le
\left(\frac{\|u\|^{p}+\|v\|^{p}}{2}\right)^{\frac{q}{p}}.
\] \end{enumerate} \end{lemma}

The scalar inequality used below is \cite[Theorem 2]{Clarkson1936}.

\begin{proof} For any $u,v\in H$, the parallelogram law implies \[
\left\|\frac{u+v}{2}\right\|^{2}
+
\left\|\frac{u-v}{2}\right\|^{2}
=
\frac{\|u\|^{2}+\|v\|^{2}}{2}.
\] Define \[
A=\left\|\frac{u+v}{2}\right\|^{2},
\qquad
B=\left\|\frac{u-v}{2}\right\|^{2},
\] so that \[
A+B=\frac{\|u\|^{2}+\|v\|^{2}}{2}.
\]

\medskip First prove (a). Since $p\ge2$, we have $\displaystyle\frac{p}{2}\ge1$. For $A,B\ge0$, this gives \[
A^{\frac{p}{2}}+B^{\frac{p}{2}}
\le
(A+B)^{\frac{p}{2}}.
\] On the other hand, convexity of the function $t\mapsto t^{\frac{p}{2}}$ on $[0,\infty)$ yields \[
(A+B)^{\frac{p}{2}}
=
\left(\frac{\|u\|^{2}+\|v\|^{2}}{2}\right)^{\frac{p}{2}}
\le
\frac{\left(\|u\|^{2}\right)^{\frac{p}{2}}
+
\left(\|v\|^{2}\right)^{\frac{p}{2}}}{2}
=
\frac{\|u\|^{p}+\|v\|^{p}}{2}.
\] Substituting $A$ and $B$ gives the result.

\medskip Finally, prove (b). Let $x=\|u\|$ and $y=\|v\|$. By the parallelogram law, the sum of our constants is \[
A+B = \frac{x^{2}+y^{2}}{2}.
\] The triangle inequality gives the following upper and lower bounds for $A$: \[
\left|\frac{x-y}{2}\right|
=
\frac{|\|u\|-\|v\||}{2}
\le
\left\|\frac{u+v}{2}\right\|
\le
\frac{\|u\|+\|v\|}{2}
=
\frac{x+y}{2}.
\] Squaring shows that $A$ belongs to the closed interval \[
I = \left[ \left(\frac{x-y}{2}\right)^{2}, \left(\frac{x+y}{2}\right)^{2} \right].
\] To bound the quantity $A^{\frac{q}{2}}+B^{\frac{q}{2}}$, consider the real-variable function $f\colon I \longrightarrow \mathbb{R}$ defined by \[
f(t) = t^{\frac{q}{2}} + \left(\frac{x^{2}+y^{2}}{2} - t\right)^{\frac{q}{2}}.
\] Since $1<p\le 2$, the conjugate exponent satisfies $q\ge 2$, so
$r=\displaystyle\frac{q}{2}\ge 1$. The function $s\mapsto s^r$ is continuous and convex on
$[0,\infty)$. Moreover, $I\subseteq\left[0,\displaystyle\frac{x^2+y^2}{2}\right]$, so the affine
maps $t\mapsto t$ and $t\mapsto \displaystyle\frac{x^2+y^2}{2-t}$ take nonnegative values
on $I$. Composing with these maps and adding shows that $f$ is
continuous and convex on all of $I$.

The maximum of a convex function on a closed interval is always attained at at least one endpoint. Evaluating $f$ at the endpoints of $I$ gives the same value in both cases: \[
f\left( \left(\frac{x+y}{2}\right)^{2} \right)
=
f\left( \left(\frac{x-y}{2}\right)^{2} \right)
=
\left(\frac{x+y}{2}\right)^{q} + \left|\frac{x-y}{2}\right|^{q}.
\] Since the constant $A \in I$, evaluation at $t=A$ yields \[
A^{\frac{q}{2}}+B^{\frac{q}{2}} = f(A) \le \left(\frac{x+y}{2}\right)^{q} + \left|\frac{x-y}{2}\right|^{q}.
\] This step geometrically reduces the vector problem to one involving nonnegative real scalars. Applying the cited scalar inequality to the real numbers $x,y\ge0$ and $1<p\le2$ gives \[
\left(\frac{x+y}{2}\right)^{q} + \left|\frac{x-y}{2}\right|^{q}
\le
\left(\frac{x^{p}+y^{p}}{2}\right)^{\frac{q}{p}}.
\] Combining the inequalities yields the desired result. \end{proof}

\begin{proposition}\label{prop: Lp(E) es uniformemente convexo} Let $(M,\mathbf{g})$ be a Riemannian manifold without boundary, and let $\mathbf{E}\to M$ be a smooth vector bundle with bundle metric $\mathbf{h}_{\mathbf{E}}$. If $1<p<\infty$, then $L^{p}(M,\mathbf{E})$ is uniformly convex. \end{proposition}

\begin{proof} Let $\mathbf{u},\mathbf{v}\in L^{p}(M,\mathbf{E})$ be such that \[
\|\mathbf{u}\|_{L^{p}(M,\mathbf{E})}=\|\mathbf{v}\|_{L^{p}(M,\mathbf{E})}=1,
\qquad
\|\mathbf{u}-\mathbf{v}\|_{L^{p}(M,\mathbf{E})}\ge \varepsilon,
\] with $\varepsilon>0$. For each $x\in M$, the fiber $(\mathbf{E}_x,\mathbf{h}_{\mathbf{E}}(x))$ is a Hilbert space, so Lemma~\ref{lem:clarkson-hilbert} applies pointwise, giving two cases: \medskip

\noindent\textit{Case $2\le p<\infty$.} Applying Lemma~\ref{lem:clarkson-hilbert}\,(a) to $\mathbf{u}(x),\mathbf{v}(x)\in \mathbf{E}_x$ and then integrating with respect to $\lambda_{\mathbf{g}}$ gives \[
\int_M \left|\frac{\mathbf{u}+\mathbf{v}}{2}\right|_{\mathbf{h}_{\mathbf{E}}}^{p}\,d\lambda_{\mathbf{g}}
+
\int_M \left|\frac{\mathbf{u}-\mathbf{v}}{2}\right|_{\mathbf{h}_{\mathbf{E}}}^{p}\,d\lambda_{\mathbf{g}}
\le
\frac{1}{2}\int_M \bigl(|\mathbf{u}|_{\mathbf{h}_{\mathbf{E}}}^{p}+|\mathbf{v}|_{\mathbf{h}_{\mathbf{E}}}^{p}\bigr)\,d\lambda_{\mathbf{g}}.
\] By the definition of the norm on $L^{p}(M,\mathbf{E})$, this is equivalent to \[
\left\|\frac{\mathbf{u}+\mathbf{v}}{2}\right\|_{L^{p}(M,\mathbf{E})}^{p}
+
\left\|\frac{\mathbf{u}-\mathbf{v}}{2}\right\|_{L^{p}(M,\mathbf{E})}^{p}
\le
\frac{\|\mathbf{u}\|_{L^{p}(M,\mathbf{E})}^{p}+\|\mathbf{v}\|_{L^{p}(M,\mathbf{E})}^{p}}{2}
=
1.
\] Moreover, \[
\left\|\frac{\mathbf{u}-\mathbf{v}}{2}\right\|_{L^{p}(M,\mathbf{E})}^{p}
=
\frac{1}{2^{p}}\|\mathbf{u}-\mathbf{v}\|_{L^{p}(M,\mathbf{E})}^{p}
\ge
\left(\frac{\varepsilon}{2}\right)^{p}.
\] Substituting into the preceding inequality, \[
\left\|\frac{\mathbf{u}+\mathbf{v}}{2}\right\|_{L^{p}(M,\mathbf{E})}^{p}
\le
1-\left(\frac{\varepsilon}{2}\right)^{p}.
\] Taking $p$th roots, \[
\left\|\frac{\mathbf{u}+\mathbf{v}}{2}\right\|_{L^{p}(M,\mathbf{E})}
\le
\left(1-\left(\frac{\varepsilon}{2}\right)^{p}\right)^{\frac{1}{p}}
\] Defining \[
\delta(\varepsilon)
:=
1-\left(1-\left(\frac{\varepsilon}{2}\right)^{p}\right)^{\frac{1}{p}}>0,
\] gives the result.

\medskip

\noindent\textit{Case $1<p<2$.} Let $q$ be the conjugate exponent of $p$. Applying Lemma~\ref{lem:clarkson-hilbert}\,(b) to $\mathbf{u}(x),\mathbf{v}(x)\in \mathbf{E}_x$ gives the pointwise inequality \[
\left|\frac{\mathbf{u}+\mathbf{v}}{2}\right|_{\mathbf{h}_{\mathbf{E}}}^{q}
+
\left|\frac{\mathbf{u}-\mathbf{v}}{2}\right|_{\mathbf{h}_{\mathbf{E}}}^{q}
\le
\left(\frac{|\mathbf{u}|_{\mathbf{h}_{\mathbf{E}}}^{p}+|\mathbf{v}|_{\mathbf{h}_{\mathbf{E}}}^{p}}{2}\right)^{\frac{q}{p}}.
\] Raising both sides to the power $\frac{p}{q}$ and integrating with respect to $\lambda_{\mathbf{g}}$ yields \[
\int_M \left( \left|\frac{\mathbf{u}+\mathbf{v}}{2}\right|_{\mathbf{h}_{\mathbf{E}}}^{q} + \left|\frac{\mathbf{u}-\mathbf{v}}{2}\right|_{\mathbf{h}_{\mathbf{E}}}^{q} \right)^{\frac{p}{q}}\,d\lambda_{\mathbf{g}}
\le
\frac{1}{2}\int_M \bigl(|\mathbf{u}|_{\mathbf{h}_{\mathbf{E}}}^{p}+|\mathbf{v}|_{\mathbf{h}_{\mathbf{E}}}^{p}\bigr)\,d\lambda_{\mathbf{g}}.
\] Since $\|\mathbf{u}\|_{L^{p}(M,\mathbf{E})}=\|\mathbf{v}\|_{L^{p}(M,\mathbf{E})}=1$, the right-hand side is exactly $1$. On the other hand, since $p<2$, the exponent satisfies $\displaystyle\frac{p}{q} = p-1 < 1$. Applying the reverse Minkowski inequality of Proposition~\ref{prop:minkowski-inversa-exponentes-menores-uno} with $r=\displaystyle\frac{p}{q}$ bounds the left-hand side from below, giving: \[
\left( \int_M \left|\frac{\mathbf{u}+\mathbf{v}}{2}\right|_{\mathbf{h}_{\mathbf{E}}}^{p}\,d\lambda_{\mathbf{g}} \right)^{\frac{q}{p}}
+
\left( \int_M \left|\frac{\mathbf{u}-\mathbf{v}}{2}\right|_{\mathbf{h}_{\mathbf{E}}}^{p}\,d\lambda_{\mathbf{g}} \right)^{\frac{q}{p}}
\le
\left( \int_M \left( \left|\frac{\mathbf{u}+\mathbf{v}}{2}\right|_{\mathbf{h}_{\mathbf{E}}}^{q} + \left|\frac{\mathbf{u}-\mathbf{v}}{2}\right|_{\mathbf{h}_{\mathbf{E}}}^{q} \right)^{\frac{p}{q}}\,d\lambda_{\mathbf{g}} \right)^{\frac{q}{p}}.
\] Since the integral on the right is bounded by $1$ and $1^{\frac{q}{p}} = 1$, in terms of the norm on $L^p(M,\mathbf{E})$ this reads \[
\left\|\frac{\mathbf{u}+\mathbf{v}}{2}\right\|_{L^{p}(M,\mathbf{E})}^{q}
+
\left\|\frac{\mathbf{u}-\mathbf{v}}{2}\right\|_{L^{p}(M,\mathbf{E})}^{q}
\le
1.
\] As in the preceding case, we have \[
\left\|\frac{\mathbf{u}-\mathbf{v}}{2}\right\|_{L^{p}(M,\mathbf{E})}^{q}
\ge
\left(\frac{\varepsilon}{2}\right)^{q},
\] which implies \[
\left\|\frac{\mathbf{u}+\mathbf{v}}{2}\right\|_{L^{p}(M,\mathbf{E})}^{q}
\le
1-\left(\frac{\varepsilon}{2}\right)^{q}.
\] Taking $q$th roots, we conclude \[
\left\|\frac{\mathbf{u}+\mathbf{v}}{2}\right\|_{L^{p}(M,\mathbf{E})}
\le
\left(1-\left(\frac{\varepsilon}{2}\right)^{q}\right)^{\frac{1}{q}}
=1-\delta(\varepsilon),
\] with the explicit choice \[
\delta(\varepsilon)
:=
1-\left(1-\left(\frac{\varepsilon}{2}\right)^{q}\right)^{\frac{1}{q}}>0.
\]

\medskip

In both cases, we have found $\delta(\varepsilon)>0$ depending only on $\varepsilon$ and $p$. This proves that $L^{p}(M,\mathbf{E})$ is uniformly convex. \end{proof} We therefore have the following: \begin{corollary}[Uniform convexity and reflexivity of $W^{m,p}(M,\mathbf{E})$]\label{Wmp uniformemente convexo}\index{uniform convexity and reflexivity} Let $(M,\mathbf{g})$ be a Riemannian manifold without boundary, and let $\mathbf{E}\to M$ be a smooth vector bundle equipped with a bundle metric $\mathbf{h}_{\mathbf{E}}$ and a connection $\nabla^{\mathbf{E}}$. If $1<p<\infty$, the Sobolev space $W^{m,p}(M,\mathbf{E})$ is uniformly convex. In particular, by Theorem~\ref{teorema milman pettis}, $W^{m,p}(M,\mathbf{E})$ is reflexive. \end{corollary}

Stability of uniform convexity under $\ell^p(\{0,\ldots,m\})$ products is \cite[Theorem 1]{Clarkson1936}. \begin{proof} By Proposition~\ref{prop: Lp(E) es uniformemente convexo}, for each $j\in\{0,\dots,m\}$, $L^{p}(M,T^{(0,j)}(TM)\otimes \mathbf{E})$ is uniformly convex. Thus the product $\displaystyle\prod_{j=0}^{m}L^{p}(M,T^{(0,j)}(TM)\otimes \mathbf{E})$ equipped with the norm $\ell^{p}(\{0,\ldots,m\})$ is uniformly convex. Uniform convexity passes to subspaces and is preserved by isometries, so Theorem~\ref{encaje isometrico en Lp sobolev haces} implies that $W^{m,p}(M,\mathbf{E})$ is uniformly convex. \end{proof} Separability will be an important property, and Lebesgue spaces inherit it from a property of measure spaces with the same name: \begin{definition}[Separable measure space]\label{def:espacio_separable}\index{separable measure space} Let $(X,\mathcal{A},\mu)$ be a measure space. We say that $(X,\mathcal{A},\mu)$ is a \emph{separable measure space} if there exists a countable family of sets $\mathcal{C} \subseteq \mathcal{A}$ such that the $\sigma$-algebra $\mathcal{A}$ is contained in the completion of the $\sigma$-algebra generated by $\mathcal{C}$ (denoted by $\overline{\sigma(\mathcal{C})}^{\,\mu}$). That is: \[
\mathcal{A} \subseteq \overline{\sigma(\mathcal{C})}^{\,\mu}.
\] In the special case where $\mathcal{A}$ is already complete, this becomes the equality: \[
\mathcal{A} = \overline{\sigma(\mathcal{C})}^{\,\mu}.
\] \end{definition} For a proof of the following lemma, see \cite{Bogachev2007}. \begin{lemma}[Separability of $L^{p}(X,\mu)$]\label{Lp separable espacio de medida separable}\index{separability} Let $(X,\mathcal{A},\mu)$ be a $\sigma$-finite measure space and let $1\le p<\infty$. If the measure space $(X,\mathcal{A},\mu)$ is separable, then $L^{p}(X,\mathcal{A},\mu)$ is separable. \end{lemma}

This tells us that $L^{p}(M)$ is separable for $1\leq p<\infty$ because $(M,\mathcal{L}_{M},\lambda_{g})$ is a separable measure space: \begin{proposition}\label{prop:medida-riemanniana-separable} Let $(M,\mathbf{g})$ be a Riemannian manifold without boundary, and let $\lambda_{\mathbf{g}}$ be the associated Riemann--Lebesgue measure. Then the measure space $(M,\mathcal L_M,\lambda_{\mathbf{g}})$ is separable and $\sigma$-finite. \end{proposition}

\begin{proof} Since $M$ is a smooth manifold, it is a second-countable topological space. In particular, its topology has a countable base of open subsets $(U_k)_{k\in\mathbb N}$ on $M$.

Let $\mathcal B$ be the Borel $\sigma$-algebra of $M$. Then \[
\mathcal B=\sigma\bigl(\{U_k\mid k\in\mathbb N\}\bigr).
\]

Moreover, by Lemma~\ref{lem:medida-riemanniana-completacion-borel-sigma-finita}, $\mathcal L_M$ agrees with the completion of $\mathcal B$ with respect to $\lambda_{\mathbf{g}}$, that is, $\mathcal L_M=\overline{\mathcal B}^{\,\lambda_{\mathbf{g}}}$, and $\lambda_{\mathbf g}$ is $\sigma$-finite.

Therefore \[
\mathcal L_M
=
\overline{\sigma(\{U_k\mid k\in\mathbb N\})}^{\,\lambda_{\mathbf{g}}},
\] showing that $(M,\mathcal L_M,\lambda_{\mathbf{g}})$ is a separable measure space. \end{proof} We use this to prove separability of $L^{p}(M,\mathbf{E})$: \begin{theorem}\label{teo: Lp(E) es separable}\index{lp e is separable@lp(E) is separable} Let $(M,\mathbf{g})$ be a Riemannian manifold without boundary, and let $\mathbf{E}\to M$ be a smooth vector bundle with bundle metric $\mathbf{h}_{\mathbf{E}}$. If $1\leq p<\infty$, then $L^{p}(M,\mathbf{E})$ is separable. \end{theorem} \begin{proof} Let $(U_{k},\phi_{k})_{k\in\mathbb N}$ be a countable atlas of $M$. Refining it by bundle trivialization domains, we may assume that for every $k$ there is a local $\mathbf{h}_{\mathbf{E}}$-orthonormal frame $(\mathbf{e}_{1,k},\dots,\mathbf{e}_{r,k})$ of $\mathbf{E}$ over $U_k$. Let $(\psi_k)_{k\in\mathbb N}$ be a smooth partition of unity subordinate to $(U_k)_{k\in\mathbb N}$.

For each $k\in\mathbb N$, since $L^{p}(U_k)$ is separable, fix a countable dense subset $\mathcal D_k\subseteq L^{p}(U_k)$. Define \[
\mathcal{D}:=
\left\{
\displaystyle\sum_{k=1}^{N}\psi_k
\left(\displaystyle\sum_{j=1}^{r} f_{j,k}\,\widetilde{\mathbf{e}}_{j,k}\right)
\ \middle|\
N\in\mathbb N,\ \ f_{j,k}\in\mathcal D_k
\right\},
\] where $\widetilde{\mathbf{e}}_{j,k}$ denotes the extension of $\mathbf{e}_{j,k}$ by zero outside $U_k$. Then $\mathcal D$ is countable and $\mathcal D\subseteq L^{p}(M,\mathbf{E})$.

Let us show that $\mathcal D$ is dense in $L^{p}(M,\mathbf{E})$. Let $\mathbf{s}\in L^{p}(M,\mathbf{E})$ and let $\varepsilon>0$. For $N\in\mathbb N$, define \[
\mathbf{s}_N:=\displaystyle\sum_{k=1}^{N}\psi_k \mathbf{s}.
\] Writing $w_N:=\displaystyle\sum_{k>N}\psi_k$, we have $0\le w_N\le 1$ and $w_N(x)\to 0$ for every $x\in M$ (since $\displaystyle\sum_{k=1}^{\infty}\psi_k(x)=1$). Moreover, \[
|\mathbf{s}-\mathbf{s}_N|_{\mathbf{h}_{\mathbf{E}}}=|w_N \mathbf{s}|_{\mathbf{h}_{\mathbf{E}}}=w_N\,|\mathbf{s}|_{\mathbf{h}_{\mathbf{E}}}\le |\mathbf{s}|_{\mathbf{h}_{\mathbf{E}}}.
\] Therefore \[
\|\mathbf{s}-\mathbf{s}_N\|_{L^{p}(M,\mathbf{E})}^{p}
=
\int_M |\mathbf{s}-\mathbf{s}_N|_{\mathbf{h}_{\mathbf{E}}}^{p}\,d\lambda_{\mathbf{g}}
=
\int_M w_N^{p}\,|\mathbf{s}|_{\mathbf{h}_{\mathbf{E}}}^{p}\,d\lambda_{\mathbf{g}}
\longrightarrow 0,
\] by Theorem~\ref{convergencia dominada}, since $w_N^{p}|\mathbf{s}|_{\mathbf{h}_{\mathbf{E}}}^{p}\to 0$ pointwise and \[
0\le w_N^{p}|\mathbf{s}|_{\mathbf{h}_{\mathbf{E}}}^{p}\le |\mathbf{s}|_{\mathbf{h}_{\mathbf{E}}}^{p}\in L^{1}(M).
\] In particular, there exists $N\in\mathbb N$ such that \[
\|\mathbf{s}-\mathbf{s}_N\|_{L^{p}(M,\mathbf{E})}<\frac{\varepsilon}{2}.
\]

Now fix $k\in\{1,\dots,N\}$. On $U_k$, write \[
\mathbf{s}\restriction_{U_k}=\displaystyle\sum_{j=1}^{r} a_{j,k}\,\mathbf{e}_{j,k},
\qquad a_{j,k}\ \text{measurable on }U_k.
\] Since the frame is $\mathbf{h}_{\mathbf{E}}$-orthonormal, for almost every $x\in U_k$ we have \[
|\mathbf{s}(x)|_{\mathbf{h}_{\mathbf{E}}}^{2}
=
\displaystyle\sum_{j=1}^{r} |a_{j,k}(x)|^{2},
\] and in particular $|a_{j,k}(x)|\le |\mathbf{s}(x)|_{\mathbf{h}_{\mathbf{E}}}$. Hence \[
\|a_{j,k}\|_{L^{p}(U_k)}^{p}
=
\int_{U_k} |a_{j,k}|^{p}\,d\lambda_{\mathbf{g}}
\le
\int_{U_k} |\mathbf{s}|_{\mathbf{h}_{\mathbf{E}}}^{p}\,d\lambda_{\mathbf{g}}
\le
\int_{M} |\mathbf{s}|_{\mathbf{h}_{\mathbf{E}}}^{p}\,d\lambda_{\mathbf{g}}
=
\|\mathbf{s}\|_{L^{p}(M,\mathbf{E})}^{p}<\infty,
\], and therefore $a_{j,k}\in L^{p}(U_k)$ for every $j\in\{1,\dots,r\}$.

Since $\mathcal D_k$ is dense in $L^{p}(U_k)$, choose $f_{j,k}\in \mathcal D_k$ such that \[
\|a_{j,k}-f_{j,k}\|_{L^{p}(U_k)}
<
\frac{\varepsilon}{2Nr},
\qquad j\in\{1,\dots,r\}.
\] Define \[
u:=\displaystyle\sum_{k=1}^{N}\psi_k
\left(\displaystyle\sum_{j=1}^{r} f_{j,k}\,\widetilde{\mathbf{e}}_{j,k}\right)\in \mathcal D.
\] For each $k\in\{1,\dots,N\}$, using $\psi_k\le 1$ and orthonormality of the frame, for almost every $x\in U_k$ we have \[
\left|
\psi_k(x)\displaystyle\sum_{j=1}^{r} (a_{j,k}(x)-f_{j,k}(x))\,\mathbf{e}_{j,k}(x)
\right|_{\mathbf{h}_{\mathbf{E}}}^{p}
=
|\psi_k(x)|^{p}
\left(\displaystyle\sum_{j=1}^{r} |a_{j,k}(x)-f_{j,k}(x)|^{2}\right)^{\frac{p}{2}}
\le
\left(\displaystyle\sum_{j=1}^{r} |a_{j,k}(x)-f_{j,k}(x)|\right)^{p}.
\] Applying the inequality \[
\left(\displaystyle\sum_{j=1}^{r} b_j\right)^{p}\le r^{p-1}\displaystyle\sum_{j=1}^{r} b_j^{p},
\qquad b_j\ge 0,
\] gives \[
\left|
\psi_k\displaystyle\sum_{j=1}^{r} (a_{j,k}-f_{j,k})\,\mathbf{e}_{j,k}
\right|_{\mathbf{h}_{\mathbf{E}}}^{p}
\le
r^{p-1}\displaystyle\sum_{j=1}^{r} |a_{j,k}-f_{j,k}|^{p}
\quad \text{in }U_k.
\] Integrating and using $\supp(\psi_k)\subset U_k$ yields \[
\left\|\psi_k\displaystyle\sum_{j=1}^{r} (a_{j,k}-f_{j,k})\,\widetilde{\mathbf{e}}_{j,k}\right\|_{L^{p}(M,\mathbf{E})}^{p}
\le
r^{p-1}\displaystyle\sum_{j=1}^{r} \|a_{j,k}-f_{j,k}\|_{L^{p}(U_k)}^{p}.
\] In particular, \[
\left\|\psi_k\displaystyle\sum_{j=1}^{r} (a_{j,k}-f_{j,k})\,\widetilde{\mathbf{e}}_{j,k}\right\|_{L^{p}(M,\mathbf{E})}
\le
r^{1-\frac{1}{p}}
\left(\displaystyle\sum_{j=1}^{r} \|a_{j,k}-f_{j,k}\|_{L^{p}(U_k)}^{p}\right)^{\frac{1}{p}}
<
\frac{\varepsilon}{2N}.
\] Summing over $k\in\{1,\dots,N\}$, \[
\|\mathbf{s}_N-u\|_{L^{p}(M,\mathbf{E})}
\le
\displaystyle\sum_{k=1}^{N}
\left\|\psi_k\displaystyle\sum_{j=1}^{r} (a_{j,k}-f_{j,k})\,\widetilde{\mathbf{e}}_{j,k}\right\|_{L^{p}(M,\mathbf{E})}
<
\displaystyle\sum_{k=1}^{N}\frac{\varepsilon}{2N}
=
\frac{\varepsilon}{2}.
\] Finally, \[
\|\mathbf{s}-u\|_{L^{p}(M,\mathbf{E})}\le \|\mathbf{s}-\mathbf{s}_N\|_{L^{p}(M,\mathbf{E})}+\|\mathbf{s}_N-u\|_{L^{p}(M,\mathbf{E})}<\varepsilon.
\] Since $\mathbf{s}$ and $\varepsilon$ were arbitrary, $\mathcal D$ is dense in $L^{p}(M,\mathbf{E})$. Consequently, $L^{p}(M,\mathbf{E})$ is separable. \end{proof}

\begin{corollary}[Separability of $W^{m,p}(M,\mathbf{E})$]\label{Wmp separable}\index{separability} Let $(M,\mathbf{g})$ be a Riemannian manifold without boundary, and let $\mathbf{E}\to M$ be a smooth vector bundle equipped with a bundle metric $\mathbf{h}_{\mathbf{E}}$ and a connection $\nabla^{\mathbf{E}}$. If $1\le p<\infty$, the Sobolev space $W^{m,p}(M,\mathbf{E})$ is separable. \end{corollary} \begin{proof} By Theorem~\ref{encaje isometrico en Lp sobolev haces}, the operator \[
J_{m}\colon W^{m,p}(M,\mathbf{E})\longrightarrow
\displaystyle \prod_{j=0}^{m}
L^{p}\bigl(M,T^{(0,j)}(TM)\otimes \mathbf{E}\bigr),
\qquad
J_m(u)=\bigl((\nabla^{\mathbf{E}})^{j}u\bigr)_{j=0}^{m},
\] is an isometry onto its image, and $J_m\bigl(W^{m,p}(M,\mathbf{E})\bigr)$ is a closed subspace.

By Theorem~\ref{teo: Lp(E) es separable}, for each $j\in\{0,\dots,m\}$ the space $L^{p}\bigl(M,T^{(0,j)}(TM)\otimes \mathbf{E}\bigr)$ is separable. Consequently, the finite product $\displaystyle \prod_{j=0}^{m}
L^{p}\!\bigl(M,T^{(0,j)}(TM)\otimes \mathbf{E}\bigr)$ is separable.

Since every subspace of a separable metric space is separable, $J_m\bigl(W^{m,p}(M,\mathbf{E})\bigr)$ is separable. Finally, since $J_m$ is an isometry onto its image (and hence a homeomorphism onto its image), we conclude that $W^{m,p}(M,\mathbf{E})$ is separable. \end{proof}

We now study the dual of $L^1(M,\mathbf{E})$. Besides extending classical $L^1(M,\mathbf E)$--$L^\infty(M,\mathbf E)$ duality to sections, this identification gives a concrete description of the weak topology of $L^1(M,\mathbf{E})$ through the bundle inner product and integration.

\begin{proposition}[Dual of $L^1(M,\mathbf{E})$]\label{dual de L1 para haces}\index{dual} Let $(M,\mathbf{g})$ be a Riemannian manifold with or without boundary, and let $\mathbf{E}\to M$ be a smooth real vector bundle of finite rank equipped with a bundle metric $\mathbf{h}_{\mathbf{E}}$. Then the map \[
J_{\mathbf{E}}\colon L^\infty(M,\mathbf{E})\longrightarrow(L^1(M,\mathbf{E}))',
\qquad
J_{\mathbf{E}}(\mathbf{v})(\mathbf{u})
:=
\int_M\langle \mathbf{u},\mathbf{v}\rangle_{\mathbf{h}_{\mathbf{E}}}\,d\lambda_{\mathbf{g}},
\] is an isometric isomorphism. In particular, $(L^1(M,\mathbf{E}))'\cong L^\infty(M,\mathbf{E})$ and $\sigma(L^1(M,\mathbf{E}),L^\infty(M,\mathbf{E}))$ agrees with $\sigma(L^1(M,\mathbf{E}),(L^1(M,\mathbf{E}))')$. \end{proposition}

\begin{proof} Riemann--Lebesgue measure is $\sigma$-finite. For $\mathbf{v}\in L^\infty(M,\mathbf{E})$, Hölder's inequality (Proposition~\ref{desigualdad de holder}) implies $|J_{\mathbf{E}}(\mathbf{v})(\mathbf{u})|\leq\|\mathbf{v}\|_{L^\infty(M,\mathbf{E})}\|\mathbf{u}\|_{L^1(M,\mathbf{E})}$ for every $\mathbf{u}\in L^1(M,\mathbf{E})$, so $J_{\mathbf{E}}$ is well defined and $\|J_{\mathbf{E}}(\mathbf{v})\|\leq\|\mathbf{v}\|_{L^\infty(M,\mathbf{E})}$.

Let us show that $J_{\mathbf{E}}$ is surjective. Let $\Lambda\in(L^1(M,\mathbf{E}))'$. Take a countable open cover $(U_j)_{j\in\mathbb N}$ of $M$, with a local trivialization of $\mathbf{E}$ defined over each $U_j$, and fix a smooth local orthonormal frame $(\mathbf{e}_{j,1},\dots,\mathbf{e}_{j,r})$ on each of these open subsets. Define the disjoint Borel sets $A_1:=U_1$ and \[
A_j
:=
U_j\setminus\bigcup_{i=1}^{j-1}U_i,
\qquad
j\geq2.
\] Then $M=\displaystyle\bigcup_{j\in\mathbb N}A_j$. For $j\in\mathbb N$ and $a\in\{1,\dots,r\}$, define \[
\Lambda_{j,a}\colon L^1(A_j)\longrightarrow\mathbb R,
\qquad
\Lambda_{j,a}(f)
:=
\Lambda(\widetilde{fe_{j,a}}),
\] where the tilde denotes extension by zero outside $A_j$. Since $\Lambda$ is continuous, \[
|\Lambda_{j,a}(f)|
=
|\Lambda(\widetilde{fe_{j,a}})|
\leq
\|\Lambda\|\,
\|\widetilde{fe_{j,a}}\|_{L^1(M,\mathbf{E})}.
\] Moreover, since the frame is orthonormal, $|\mathbf{e}_{j,a}|_{\mathbf{h}_{\mathbf{E}}}=1$ on $A_j$, so \[
\|\widetilde{fe_{j,a}}\|_{L^1(M,\mathbf{E})}
=
\int_{A_j}|fe_{j,a}|_{\mathbf{h}_{\mathbf{E}}}\,d\lambda_{\mathbf{g}}
=
\int_{A_j}|f|\,d\lambda_{\mathbf{g}}
=
\|f\|_{L^1(A_j)}.
\] Consequently, $|\Lambda_{j,a}(f)|\leq\|\Lambda\|\|f\|_{L^1(A_j)}$, and hence $\|\Lambda_{j,a}\|\leq\|\Lambda\|$. By Theorem~\ref{teo:dualidad-reflexividad-Lp}, applied to $L^1(A_j)$, there exists a unique $v_j^a\in L^\infty(A_j)$ such that \[
\Lambda_{j,a}(f)
=
\int_{A_j}fv_j^a\,d\lambda_{\mathbf{g}},
\qquad
\|v_j^a\|_{L^\infty(A_j)}
=
\|\Lambda_{j,a}\|
\leq
\|\Lambda\|.
\]

Define a measurable section $\mathbf{v}$ by $\mathbf{v}\restriction_{A_j}:=\displaystyle\sum_{a=1}^rv_j^a\mathbf{e}_{j,a}$. Since the frames are orthonormal, $|\mathbf{v}|_{\mathbf{h}_{\mathbf{E}}}^2
=
\displaystyle\sum_{a=1}^r|v_j^a|^2
\leq
r\|\Lambda\|^2$ almost everywhere on $A_j$, giving $|\mathbf{v}|_{\mathbf{h}_{\mathbf{E}}}\leq\sqrt r\,\|\Lambda\|$ almost everywhere and therefore $\mathbf{v}\in L^\infty(M,\mathbf{E})$.

For $\mathbf{u}\in L^1(M,\mathbf{E})$, write $u_j^a:=\langle \mathbf{u},\mathbf{e}_{j,a}\rangle_{\mathbf{h}_{\mathbf{E}}}$ on $A_j$ and define $\mathbf{u}_N:=\displaystyle\sum_{j=1}^N\mathbf{1}_{A_j}\mathbf{u}$. Since \[
\|\mathbf{u}-\mathbf{u}_N\|_{L^1(M,\mathbf{E})}
=
\displaystyle\sum_{j>N}\int_{A_j}|\mathbf{u}|_{\mathbf{h}_{\mathbf{E}}}\,d\lambda_{\mathbf{g}}
\longrightarrow0,
\] we have $\mathbf{u}_N\to \mathbf{u}$ in $L^1(M,\mathbf{E})$. Moreover, $\mathbf{u}_N
=
\displaystyle\sum_{j=1}^N\displaystyle\sum_{a=1}^r
\widetilde{u_j^ae_{j,a}}$.

By linearity of $\Lambda$ and the definition of $\Lambda_{j,a}$, \[
\Lambda(\mathbf{u}_N)
=
\displaystyle\sum_{j=1}^N\displaystyle\sum_{a=1}^r
\Lambda\bigl(\widetilde{u_j^ae_{j,a}}\bigr)
=
\displaystyle\sum_{j=1}^N\displaystyle\sum_{a=1}^r
\Lambda_{j,a}(u_j^a)
=
\displaystyle\sum_{j=1}^N\displaystyle\sum_{a=1}^r
\int_{A_j}u_j^av_j^a\,d\lambda_{\mathbf{g}}.
\] Since the frames are orthonormal, $\langle \mathbf{u},\mathbf{v}\rangle_{\mathbf{h}_{\mathbf{E}}}
=
\displaystyle\sum_{a=1}^ru_j^av_j^a$ on each $A_j$, so \[
\displaystyle\sum_{j=1}^N\displaystyle\sum_{a=1}^r
\int_{A_j}u_j^av_j^a\,d\lambda_{\mathbf{g}}
=
\int_M
\mathbf{1}_{\bigcup_{j=1}^NA_j}
\langle \mathbf{u},\mathbf{v}\rangle_{\mathbf{h}_{\mathbf{E}}}\,d\lambda_{\mathbf{g}}.
\] Moreover, \[
\mathbf{1}_{\bigcup_{j=1}^NA_j}
\langle \mathbf{u},\mathbf{v}\rangle_{\mathbf{h}_{\mathbf{E}}}
\longrightarrow
\langle \mathbf{u},\mathbf{v}\rangle_{\mathbf{h}_{\mathbf{E}}}
\quad\text{almost everywhere},
\] and \[
\left|
\mathbf{1}_{\bigcup_{j=1}^NA_j}
\langle \mathbf{u},\mathbf{v}\rangle_{\mathbf{h}_{\mathbf{E}}}
\right|
\leq
|\mathbf{u}|_{\mathbf{h}_{\mathbf{E}}}|\mathbf{v}|_{\mathbf{h}_{\mathbf{E}}}\in L^1(M),
\] since $\mathbf{u}\in L^1(M,\mathbf{E})$ and $\mathbf{v}\in L^\infty(M,\mathbf{E})$. By Theorem~\ref{convergencia dominada} and continuity of $\Lambda$, \[
\Lambda(\mathbf{u})
=
\lim_{N\to\infty}\Lambda(\mathbf{u}_N)
=
\lim_{N\to\infty}
\displaystyle\sum_{j=1}^N\displaystyle\sum_{a=1}^r
\int_{A_j}u_j^av_j^a\,d\lambda_{\mathbf{g}}
=
\int_M\langle \mathbf{u},\mathbf{v}\rangle_{\mathbf{h}_{\mathbf{E}}}\,d\lambda_{\mathbf{g}}.
\] Consequently, $\Lambda=J_{\mathbf{E}}(\mathbf{v})$.

Finally, prove that $J_{\mathbf{E}}$ is isometric. We already know that $\|J_{\mathbf{E}}(\mathbf{v})\|\leq\|\mathbf{v}\|_{L^\infty(M,\mathbf{E})}$. Suppose that $\mathbf{v}\neq0$ and let $0<\varepsilon<\|\mathbf{v}\|_{L^\infty(M,\mathbf{E})}$. By the definition of essential supremum, the set \[
A_\varepsilon
:=
\left\{
x\in M
\;\middle|\;
|\mathbf{v}(x)|_{\mathbf{h}_{\mathbf{E}}}>
\|\mathbf{v}\|_{L^\infty(M,\mathbf{E})}-\varepsilon
\right\}
\] has positive measure, since otherwise $|\mathbf{v}|_{\mathbf{h}_{\mathbf{E}}}\leq\|\mathbf{v}\|_{L^\infty(M,\mathbf{E})}-\varepsilon$ almost everywhere, contradicting the definition of $\|\mathbf{v}\|_{L^\infty(M,\mathbf{E})}$. Since $\lambda_{\mathbf{g}}$ is $\sigma$-finite, there exist measurable sets $(C_k)_{k\in\mathbb N}$ such that $M=\displaystyle\bigcup_{k\in\mathbb N}C_k$ and $\lambda_{\mathbf{g}}(C_k)<\infty$. Since $A_\varepsilon=\displaystyle\bigcup_{k\in\mathbb N}(A_\varepsilon\cap C_k)$ has positive measure, there exists $k\in\mathbb N$ such that $B_\varepsilon:=A_\varepsilon\cap C_k$ satisfies $0<\lambda_{\mathbf{g}}(B_\varepsilon)<\infty$.

Define \[
\mathbf{u}_\varepsilon(x)
:=
\begin{cases}
\displaystyle
\frac{\mathbf{v}(x)}
{\lambda_{\mathbf{g}}(B_\varepsilon)|\mathbf{v}(x)|_{\mathbf{h}_{\mathbf{E}}}},
&x\in B_\varepsilon,\\[6pt]
0_x,
&x\notin B_\varepsilon.
\end{cases}
\] Since $B_\varepsilon\subseteq A_\varepsilon$ and $\varepsilon<\|\mathbf{v}\|_{L^\infty(M,\mathbf{E})}$, we have $|\mathbf{v}|_{\mathbf{h}_{\mathbf{E}}}>0$ on $B_\varepsilon$. Moreover, \[
\|\mathbf{u}_\varepsilon\|_{L^1(M,\mathbf{E})}
=
\frac{1}{\lambda_{\mathbf{g}}(B_\varepsilon)}
\int_{B_\varepsilon}1\,d\lambda_{\mathbf{g}}
=
1,
\] and \[
J_{\mathbf{E}}(\mathbf{v})(\mathbf{u}_\varepsilon)
=
\frac{1}{\lambda_{\mathbf{g}}(B_\varepsilon)}
\int_{B_\varepsilon}|\mathbf{v}|_{\mathbf{h}_{\mathbf{E}}}\,d\lambda_{\mathbf{g}}
>
\|\mathbf{v}\|_{L^\infty(M,\mathbf{E})}-\varepsilon.
\] By the definition of the operator norm, \[
\|J_{\mathbf{E}}(\mathbf{v})\|
=
\sup_{\substack{\mathbf u\in L^1(M,\mathbf E)\\\|\mathbf u\|_{L^1(M,\mathbf E)}\leq1}}|J_{\mathbf{E}}(\mathbf{v})(\mathbf{u})|
\geq
|J_{\mathbf{E}}(\mathbf{v})(\mathbf{u}_\varepsilon)|
>
\|\mathbf{v}\|_{L^\infty(M,\mathbf{E})}-\varepsilon.
\] Since this inequality holds for every $0<\varepsilon<\|\mathbf{v}\|_{L^\infty(M,\mathbf{E})}$, letting $\varepsilon\to0$ gives $\|J_{\mathbf{E}}(\mathbf{v})\|\geq\|\mathbf{v}\|_{L^\infty(M,\mathbf{E})}$. Consequently, $\|J_{\mathbf{E}}(\mathbf{v})\|=\|\mathbf{v}\|_{L^\infty(M,\mathbf{E})}$, and $J_{\mathbf{E}}$ is an isometric isomorphism. \end{proof}

For a smooth \emph{real} vector bundle $\mathbf{E}\to M$ of finite rank equipped with a bundle metric $\mathbf{h}_{\mathbf{E}}$, consider on $L^1(M,\mathbf{E})\times L^\infty(M,\mathbf{E})$ the duality \[
\mathfrak b_{\mathbf{E}}(u,\mathbf{v})
:=
\int_M\langle u,\mathbf{v}\rangle_{\mathbf{h}_{\mathbf{E}}}\,d\lambda_{\mathbf{g}}.
\] By Proposition~\ref{dual de L1 para haces}, the pairing $\mathfrak b_{\mathbf{E}}$ identifies $L^\infty(M,\mathbf{E})$ isometrically with $(L^1(M,\mathbf{E}))'$. Consequently, \[
\sigma(L^1(M,\mathbf{E}),L^\infty(M,\mathbf{E}))
=
\sigma(L^1(M,\mathbf{E}),(L^1(M,\mathbf{E}))'),
\] that is, $\sigma(L^1(M,\mathbf{E}),L^\infty(M,\mathbf{E}))$ agrees with the weak topology of $L^1(M,\mathbf{E})$. Thus a net $(u_\gamma)_\gamma$ converges to $u$ in this topology if and only if \[
\int_M\langle u_\gamma,\mathbf{v}\rangle_{\mathbf{h}_{\mathbf{E}}}\,d\lambda_{\mathbf{g}}
\longrightarrow
\int_M\langle u,\mathbf{v}\rangle_{\mathbf{h}_{\mathbf{E}}}\,d\lambda_{\mathbf{g}}
\] for every $\mathbf{v}\in L^\infty(M,\mathbf{E})$.

If $\mathbf{E}$ is complex Hermitian, the complex-linear formulation replaces the second factor with $L^\infty(M,\overline{\mathbf{E}})$ and uses $\displaystyle \mathfrak b_{\mathbf{E}}(u,\overline v)=\int_M\mathbf{h}_{\mathbf{E}}(u,v)\,d\lambda_{\mathbf{g}}$; this version, for every $1\leq p<\infty$, is proved in Corollary~\ref{cor:dualidad-Lp-haz-hermitiano} of Appendix~B.

\subsection{Compatibility with restrictions and local expressions for covariant derivatives on vector bundles} Another important property of Sobolev sections is that restriction to an open subset gives a Sobolev section on that open subset. \begin{proposition}\label{prop: restriccion de seccion de Sobolev} Let $(M,\mathbf{g})$ be a Riemannian manifold without boundary, and let $\mathbf{E}\to M$ be a smooth vector bundle equipped with a bundle metric $\mathbf{h}_{\mathbf{E}}$ and a smooth connection $\nabla^{\mathbf E}$. If $U\subseteq M$ is open, $\mathbf{u}\in W^{m,p}(M,\mathbf{E})$, and $s\in\{0,\dots,m\}$, then $\mathbf{u}\restriction_U\in W^{m,p}(U,\mathbf{E}\restriction_U)$ and $\nabla_w^s(\mathbf{u}\restriction_U)
=
(\nabla_w^s\mathbf{u})\restriction_U$ almost everywhere on $U$. \end{proposition} \begin{proof} Let $s\leq m$. We must prove that $\mathbf{u}\restriction_{U}$ has a weak covariant derivative of order $s$ on $L^{p}(U,T^{(0,s)}(TU)\otimes \mathbf{E}\restriction_{U})$. To this end, take $\boldsymbol{\sigma} \in\Gamma_{c}(T^{(0,s)}(TU)\otimes \mathbf{E}\restriction_{U})$. We may extend $\boldsymbol{\sigma}$ to $\widetilde{\boldsymbol{\sigma}}\in \Gamma_{c}(T^{(0,s)}(TM)\otimes \mathbf{E})$ as in Lemma~\ref{lema: principio de localizacion para haces vectoriales}. On the other hand, since $\mathbf{u}\in W^{m,p}(M,\mathbf{E})$, there exists $\nabla_{w}^{s}\mathbf{u}$, the $s$th weak covariant derivative of $\mathbf{u}$. It satisfies \[\int_{M}\langle \nabla_{w}^{s}\mathbf{u},\widetilde{\boldsymbol{\sigma}}\rangle_{\mathbf{g},\mathbf{h}_{\mathbf{E}}}d\lambda_{\mathbf{g}}=\int_{M}\langle \mathbf{u},(\nabla^{s})^{*}\widetilde{\boldsymbol{\sigma}}\rangle_{\mathbf{h}_{\mathbf{E}}}d\lambda_{\mathbf{g}}.\] Since $\boldsymbol{\sigma} \equiv \widetilde{\boldsymbol{\sigma}}$ on $U$, which is open, Remark~\ref{remark: pdo es local implica igualdad de pdos si hay igualdad de funciones} gives $(\nabla^{s})^{*}\widetilde{\boldsymbol{\sigma} }\equiv (\nabla^{s})^{*}\boldsymbol{\sigma}$ on $U$. Then \[\int_{U}\langle (\nabla_{w}^{s}\mathbf{u})\restriction_{U},\boldsymbol{\sigma}\rangle_{\mathbf{g},\mathbf{h}_{\mathbf{E}}}d\lambda_{\mathbf{g}}=\int_{M}\langle \nabla_{w}^{s}\mathbf{u},\widetilde{\boldsymbol{\sigma}}\rangle_{\mathbf{g},\mathbf{h}_{\mathbf{E}}}d\lambda_{\mathbf{g}}=\int_{M}\langle \mathbf{u},(\nabla^{s})^{*}\widetilde{\boldsymbol{\sigma}}\rangle_{\mathbf{h}_{\mathbf{E}}}d\lambda_{\mathbf{g}}=\int_{U}\langle \mathbf{u}\restriction_{U},(\nabla^{s})^{*}\boldsymbol{\sigma} \rangle_{\mathbf{h}_{\mathbf{E}}}d\lambda_{\mathbf{g}},\] where the last equality holds because Proposition~\ref{pdo es local} gives $\supp((\nabla^{s})^{*}\widetilde{\boldsymbol{\sigma} })\subseteq\supp(\widetilde{\boldsymbol{\sigma}})=\supp(\boldsymbol{\sigma})\subseteq U$ and $(\nabla^{s})^{*}\widetilde{\boldsymbol{\sigma} }\equiv (\nabla^{s})^{*}\boldsymbol{\sigma}$ on $U$. Thus $(\nabla_w^s\mathbf{u})\restriction_U$ satisfies the identity characterizing the $s$th weak covariant derivative of $\mathbf{u}\restriction_U$. By uniqueness of the weak covariant derivative, $\nabla_w^s(\mathbf{u}\restriction_U)
=
(\nabla_w^s\mathbf{u})\restriction_U$ almost everywhere on $U$. Since $s\le m$ was arbitrary, $\mathbf{u}\restriction_U\in W^{m,p}(U,\mathbf{E}\restriction_U)$. \end{proof} Studying local properties of Sobolev spaces is important because it gives a way to work with them computationally. To achieve this, the reader should first become thoroughly familiar with coordinate expressions for higher-order covariant derivatives on smooth vector bundles. Such expressions are common in the literature in certain cases, but higher orders are less often treated because of their complexity. In this section, we calculate these local expressions in detail.

\begin{lemma}\label{lema: derivada covariante en coordenadas haz tensorial y E} Let $(M,\mathbf{g})$ be a Riemannian manifold with or without boundary with a connection $\nabla$ on $TM$, and let $\mathbf{E}$ be a smooth vector bundle with a connection $\nabla^{\mathbf{E}}$. Also denote by $\nabla$ the induced connection on the tensor bundle $T^{(k,l)}(TM)\otimes \mathbf{E}$ characterized by Proposition~\ref{conexion en tensores, duales y endomorfismos}. Let $U\subseteq M$ be an open subset on which $(\mathbf{E}_{1},\dots,\mathbf{E}_{n})$ is a local frame of $TM$ over $U$, $(\boldsymbol{\varepsilon}^{1},\dots,\boldsymbol{\varepsilon}^{n})$ its associated dual coframe, and $(\mathbf{e}_{1},\dots,\mathbf{e}_{r})$ a local frame of $\mathbf{E}$ over $U$. If $\mathbf{F}\in \Gamma(T^{(k,l)}(TM)\otimes \mathbf{E})$ is expressed locally as \[
\mathbf{F} = F^{i_{1}\dots i_{k}a}_{j_{1}\dots j_{l}} \, \mathbf{E}_{i_{1}}\otimes \cdots \otimes \mathbf{E}_{i_{k}}\otimes \boldsymbol{\varepsilon}^{j_{1}}\otimes \cdots \otimes \boldsymbol{\varepsilon}^{j_{l}}\otimes \mathbf{e}_{a},
\] then the components of $\nabla \mathbf{F}$ in these frames are \[(\nabla_{\mathbf{E}_c}\mathbf{F})^{i_{1}\dots i_{k}a}_{j_{1}\dots j_{l}}
= \mathbf{E}_{c}\big(F^{i_{1}\dots i_{k}a}_{j_{1}\dots j_{l}}\big)
+ \displaystyle\sum_{t=1}^{k}\Gamma^{i_{t}}_{cs}\, F^{i_{1}\dots s\dots i_{k}a}_{j_{1}\dots j_{l}}
 - \displaystyle\sum_{t=1}^{l}\Gamma^{s}_{c j_{t}}\, F^{i_{1}\dots i_{k}\,a}_{j_{1}\dots s\dots j_{l}}
+ \displaystyle\sum_{b=1}^{r}(A_{c})^{a}_{b}\, F^{i_{1}\dots i_{k}b}_{j_{1}\dots j_{l}},\] where the connection coefficients are given by $\nabla_{\mathbf{E}_{c}}\mathbf{E}_{j}=\Gamma^{s}_{c j}\mathbf{E}_{s}$ and $\nabla^{\mathbf{E}}_{\mathbf{E}_{c}}\mathbf{e}_{b}=(A_{c})^{a}_{b}\,\mathbf{e}_{a}$. \end{lemma} \begin{proof} The covariant derivative of $\mathbf{F}$ in the direction $\mathbf{E}_{c}$ is \[
\nabla_{\mathbf{E}_c} \mathbf{F} = \nabla_{\mathbf{E}_c} \left( \displaystyle\sum_{a=1}^r \mathbf{F}^{(a)} \otimes \mathbf{e}_a \right),
\] where $\mathbf{F}^{(a)} \in \Gamma(T^{(k,l)}(TM))$ is the tensor field with components $(\mathbf{F}^{(a)})^{i_1 \dots i_k}_{j_1 \dots j_l} = F^{i_1 \dots i_k a}_{j_1 \dots j_l}$. By Proposition~\ref{conexion en tensores, duales y endomorfismos}, applying the Leibniz rule gives \begin{equation} \label{eq:leibniz-expansion}
\nabla_{\mathbf{E}_c} \mathbf{F} = \displaystyle\sum_{a=1}^r \left( (\nabla_{\mathbf{E}_c} \mathbf{F}^{(a)}) \otimes \mathbf{e}_a + \mathbf{F}^{(a)} \otimes \nabla^{\mathbf{E}}_{\mathbf{E}_c} \mathbf{e}_a \right).
\end{equation} By Lemma~\ref{lema:formula-derivada-covariante-coords}, the coefficients of the first term in \eqref{eq:leibniz-expansion} are \[
(\nabla_{\mathbf{E}_c} \mathbf{F}^{(a)})^{i_1 \dots i_k}_{j_1 \dots j_l} = \mathbf{E}_c(F^{i_1 \dots i_k a}_{j_1 \dots j_l}) + \displaystyle\sum_{t=1}^{k}\Gamma^{i_{t}}_{cs}\, F^{i_{1}\dots s\dots i_{k}a}_{j_{1}\dots j_{l}} - \displaystyle\sum_{t=1}^{l}\Gamma^{s}_{c j_{t}}\, F^{i_{1}\dots i_{k}\,a}_{j_{1}\dots s\dots j_{l}}.
\] For the second term in \eqref{eq:leibniz-expansion}, substitute $\nabla^{\mathbf{E}}_{\mathbf{E}_c} \mathbf{e}_b = (A_c)^a_b \mathbf{e}_a$. Changing the summation index for the local frame $(\mathbf{e}_{1},\dots,\mathbf{e}_{r})$ to $b$ in \eqref{eq:leibniz-expansion} yields \[
\displaystyle\sum_{b=1}^r \mathbf{F}^{(b)} \otimes \nabla^{\mathbf{E}}_{\mathbf{E}_c} \mathbf{e}_b = \displaystyle\sum_{b=1}^r \mathbf{F}^{(b)} \otimes \left( \displaystyle\sum_{a=1}^r (A_c)^a_b \mathbf{e}_a \right) = \displaystyle\sum_{a=1}^r \left( \displaystyle\sum_{b=1}^r (A_c)^a_b \mathbf{F}^{(b)} \right) \otimes \mathbf{e}_a.
\] The component of this term corresponding to $\mathbf{e}_{a}$ is $\displaystyle\sum_{b=1}^r (A_c)^a_b F^{i_1 \dots i_k b}_{j_1 \dots j_l}$, completing the proof. \end{proof} This gives the following local expression for higher-order derivatives: \begin{lemma}[Local expression for the $s$th covariant derivative] \label{lem:local-expression-higher-order}\index{covariant derivative!local expression for the \(s\)th derivative} Let $M$ be a smooth manifold with or without boundary equipped with a connection $\nabla$ on $TM$, and let $\pi_{\mathbf{E}}\colon \mathbf{E}\longrightarrow M$ be a smooth vector bundle of rank $r$ equipped with a connection $\nabla^{\mathbf{E}}$. Let $(U,x^{1},\dots,x^{n})$ be a smooth chart of $M$ and $(\mathbf{e}_{1},\dots,\mathbf{e}_{r})$ a local frame of $\mathbf{E}$ over $U$. Denote the connection coefficients by \[
\nabla_{\boldsymbol{\partial}_{i}}\boldsymbol{\partial}_{j} = \displaystyle\sum_{k=1}^{n}\Gamma^{k}_{ij}\,\boldsymbol{\partial}_{k},
\qquad
\nabla^{\mathbf{E}}_{\boldsymbol{\partial}_{i}}\mathbf{e}_{b} = \displaystyle\sum_{a=1}^{r} (A_{i})^{a}_{b}\,\mathbf{e}_{a}.
\] If $\mathbf{u}\in \Gamma(\mathbf{E})$ has coordinate expression $\mathbf{u} = \displaystyle\sum_{b=1}^{r}u^{b}\mathbf{e}_{b}$ on $U$, then for every integer $s \ge 1$, the components of the $s$th covariant derivative $\nabla^{s}\mathbf{u} \in \Gamma(T^{(0,s)}(TM)\otimes \mathbf{E})$ are given by \begin{equation} \label{eq:higher-order-formula}
(\nabla^{s}\mathbf{u})^{a}_{i_{1}\dots i_{s}}
=
\partial_{i_{1}}\cdots\partial_{i_{s}}u^{a}
+
\displaystyle\sum_{b=1}^{r} \displaystyle\sum_{|\beta|\le s-1}
\bigl(B_{\beta}\bigr)^{a}_{b\,i_{1}\dots i_{s}}
\,\partial^{\beta}u^{b},
\end{equation} where the coefficients $\bigl(B_{\beta}\bigr)^{a}_{b\,i_{1}\dots i_{s}}$ are smooth functions on $U$ depending only on the connection coefficients $\Gamma^{k}_{ij}$, $(A_{i})^{a}_{b}$, and their partial derivatives through order $s-1$. \end{lemma} \begin{proof} We prove the assertion by induction on $s$. For $s=1$, the expression for the connection on $\mathbf{E}$ in the chosen frame gives \[
(\nabla \mathbf{u})^a_i
=
\partial_i u^a+
\sum_{b=1}^{r}(A_i)^a{}_b u^b.
\] This is \eqref{eq:higher-order-formula} with $\bigl(B^{(1)}_0\bigr)^a_{b,i}=(A_i)^a{}_b$.

Now suppose that, for some $s\geq1$, there exist smooth coefficients $B^{(s)}_\beta$, with $|\beta|\leq s-1$, for which \begin{equation}\label{eq:hipotesis-inductiva-derivada-covariante-local}
(\nabla^s\mathbf{u})^a_{i_2\dots i_{s+1}}
=
\partial_{i_2}\cdots\partial_{i_{s+1}}u^a
+
\sum_{b=1}^{r}\sum_{|\beta|\leq s-1}
\bigl(B^{(s)}_\beta\bigr)^a_{b,i_2\dots i_{s+1}}
\partial^\beta u^b.
\end{equation} As part of the induction hypothesis, also suppose these coefficients depend only on $\Gamma^k_{ij}$, $(A_i)^a{}_b$, and their partial derivatives of order at most $s-1$.

Set $\mathbf{G}:=\nabla^s\mathbf{u}$ and abbreviate $I=(i_2,\dots,i_{s+1})$. For $t\in\{2,\dots,s+1\}$, denote by $I_t(d)$ the list obtained from $I$ by replacing $i_t$ with $d$. The formula for the induced connection on $T^{(0,s)}(TM)\otimes\mathbf{E}$ gives \begin{equation}\label{eq:paso-inductivo-derivada-covariante-local}
\begin{aligned}
(\nabla\mathbf{G})^a_{i_1I}
={}&
\partial_{i_1}G^a_I
-\sum_{t=2}^{s+1}\sum_{d=1}^{n}
\Gamma^d_{i_1i_t}G^a_{I_t(d)}
+\sum_{c=1}^{r}(A_{i_1})^a{}_cG^c_I.
\end{aligned}
\end{equation} Substitute \eqref{eq:hipotesis-inductiva-derivada-covariante-local} into each of the three terms of \eqref{eq:paso-inductivo-derivada-covariante-local}. In the first, the Leibniz rule gives \begin{align*}
\partial_{i_1}G^a_I
=
\partial_{i_1}\partial_{i_2}\cdots\partial_{i_{s+1}}u^a+
\sum_{b=1}^{r}\sum_{|\beta|\leq s-1}
\partial_{i_1}\!\left(
\bigl(B^{(s)}_\beta\bigr)^a_{b,I}
\right)\partial^\beta u^b+
\sum_{b=1}^{r}\sum_{|\beta|\leq s-1}
\bigl(B^{(s)}_\beta\bigr)^a_{b,I}
\partial^{\beta+e_{i_1}}u^b,
\end{align*} where $e_{i_1}$ is the multi-index whose only nonzero component is the $i_1$th. The second term becomes \begin{align*}
-\sum_{t=2}^{s+1}\sum_{d=1}^{n}
\Gamma^d_{i_1i_t}G^a_{I_t(d)}
=
-\sum_{t=2}^{s+1}\sum_{d=1}^{n}
\Gamma^d_{i_1i_t}
\partial_{i_2}\cdots\partial_{i_{t-1}}
\partial_d
\partial_{i_{t+1}}\cdots\partial_{i_{s+1}}u^a-
\sum_{t=2}^{s+1}\sum_{d=1}^{n}
\sum_{b=1}^{r}\sum_{|\beta|\leq s-1}
\Gamma^d_{i_1i_t}
\bigl(B^{(s)}_\beta\bigr)^a_{b,I_t(d)}
\partial^\beta u^b.
\end{align*} Finally, the term corresponding to the connection on $\mathbf{E}$ is \begin{align*}
\sum_{c=1}^{r}(A_{i_1})^a{}_cG^c_I
=
\sum_{c=1}^{r}(A_{i_1})^a{}_c
\partial_{i_2}\cdots\partial_{i_{s+1}}u^c+
\sum_{c=1}^{r}\sum_{b=1}^{r}
\sum_{|\beta|\leq s-1}
(A_{i_1})^a{}_c
\bigl(B^{(s)}_\beta\bigr)^c_{b,I}
\partial^\beta u^b.
\end{align*}

These three equalities contain exactly one term with $s+1$ derivatives of a component of $\mathbf{u}$, namely \[
\partial_{i_1}\partial_{i_2}\cdots\partial_{i_{s+1}}u^a.
\] The terms in the third line of the first expansion contain derivatives of order $|\beta|+1\leq s$; the first summand of the second expansion and the first of the third contain derivatives of order $s$; and all remaining summands contain derivatives of order at most $s-1$.

For fixed $a,b,i_1,\dots,i_{s+1}$ and a multi-index $\gamma$ with $|\gamma|\leq s$, define $\bigl(B^{(s+1)}_\gamma\bigr)^a_{b,i_1\dots i_{s+1}}$ as the sum of the coefficients multiplying $\partial^\gamma u^b$ in the preceding three expansions. This is a finite sum of products of $\Gamma^k_{ij}$, $(A_i)^a{}_b$, $B^{(s)}_\beta$, and $\partial_{i_1}B^{(s)}_\beta$. By the induction hypothesis, $B^{(s)}_\beta$ depends on the connection coefficients and their derivatives through order $s-1$, while $\partial_{i_1}B^{(s)}_\beta$ depends on them through order $s$. Consequently, $B^{(s+1)}_\gamma$ is smooth and depends only on $\Gamma^k_{ij}$, $(A_i)^a{}_b$, and their partial derivatives of order at most $s$. With this definition, we obtain \[
(\nabla^{s+1}\mathbf{u})^a_{i_1\dots i_{s+1}}
=
\partial_{i_1}\cdots\partial_{i_{s+1}}u^a
+\sum_{b=1}^{r}\sum_{|\gamma|\leq s}
\bigl(B^{(s+1)}_\gamma\bigr)^a_{b,i_1\dots i_{s+1}}
\partial^\gamma u^b,
\] which is the assertion for $s+1$ and completes the induction. \end{proof} \subsection{The integration-by-parts formula for vector bundles and the formal adjoint} One of the most important examples of a formal adjoint is that of the $s$th covariant derivative. This adjoint can be computed globally and explicitly from an integration-by-parts formula. This is natural, since the concept of formal adjoint is motivated precisely by the integration-by-parts formulas that gave rise to weak derivatives in the Euclidean case.

In this section, we prove an original result found in \cite{daniel2026geometricintegrationpartssobolev} that made possible an intrinsic, unified treatment of Sobolev spaces. Its derivation was initially motivated by \cite{chan2024meyersserrin}, which sought to establish the Meyers--Serrin theorem on manifolds in certain cases by giving explicit expressions for first- and second-order weak covariant derivatives. Subsequent analysis of these ideas culminated in the work mentioned first. The result is as follows:

\begin{theorem}[Integration-by-parts formula for higher-order covariant derivatives on vector bundles]\label{teo:green-m-esima-covariante-E}\index{integration-by-parts formula@integration-by-parts formula!for higher-order covariant derivatives} Let $(M,\mathbf{g})$ be a Riemannian manifold with or without boundary and Levi--Civita connection $\nabla$. Let $\mathbf{E}\longrightarrow M$ be a smooth vector bundle equipped with a bundle metric $\mathbf{h}_{\mathbf{E}}$ and a compatible connection $\nabla^{\mathbf{E}}$. Also denote by $\nabla$ the induced connection on the tensor bundle $T^{(k,l)}(TM)\otimes \mathbf{E}$ characterized by Proposition~\ref{conexion en tensores, duales y endomorfismos}. If $s\geq 1$, $\mathbf{F}\in \Gamma_{c}\bigl(T^{(k,l)}(TM)\otimes \mathbf{E}\bigr)$, and $\mathbf{G}\in \Gamma_{c}\bigl(T^{(k,l+s)}(TM)\otimes \mathbf{E}\bigr)$, then \[
\int_{M}\langle \nabla^{s}\mathbf{F},\mathbf{G}\rangle_{\mathbf{g},\mathbf{h}_{\mathbf{E}}}\,d\lambda_{\mathbf{g}}
=
(-1)^{s}\int_{M}\big\langle \mathbf{F},(\operatorname{tr}_{\mathbf{g}}\circ \nabla)^{s}\mathbf{G}\big\rangle_{\mathbf{g},\mathbf{h}_{\mathbf{E}}}\,d\lambda_{\mathbf{g}}
\] \[
\qquad
+\displaystyle\sum_{j=0}^{s-1}(-1)^{\,s-1-j}
\int_{\partial M}\Big\langle \nabla^{j}\mathbf{F},\,\iota_{\boldsymbol{\nu}}\bigl((\operatorname{tr}_{\mathbf{g}}\circ \nabla)^{\,s-1-j}\mathbf{G}\bigr)\Big\rangle_{\mathbf{g},\mathbf{h}_{\mathbf{E}}}\,d\lambda_{\widetilde{\mathbf{g}}},
\] where $\boldsymbol{\nu}$ is the outward unit normal to $\partial M$ and $\iota_{\boldsymbol{\nu}}$ denotes contraction of $\boldsymbol{\nu}$ in the first covariant index. More precisely, if \[
\mathbf{H}\in \Gamma\bigl(T^{(k,m+1)}(TM)\otimes \mathbf{E}\bigr),
\] then $\iota_{\boldsymbol{\nu}}\mathbf{H}\in \Gamma\bigl((T^{(k,m)}(TM)\otimes \mathbf{E})\restriction_{\partial M}\bigr)$ and, in local coordinates on the boundary, \[
\bigl(\iota_{\boldsymbol{\nu}}\mathbf{H}\bigr)^{i_{1}\dots i_{k}a}_{j_{1}\dots j_{m}}
:=
\nu^{q}\,H^{i_{1}\dots i_{k}a}_{qj_{1}\dots j_{m}},
\] where $\boldsymbol{\nu}=\nu^{q}\boldsymbol{\partial}_{q}$ in those coordinates. If \(\partial M=\varnothing\), the sum of integrals over \(\partial M\) is interpreted as zero, and no choice of normal is involved.

Here $\operatorname{tr}_{\mathbf{g}}$ denotes contraction, using $\mathbf{g}$, of the first two covariant indices of $\nabla \mathbf{H}$. More precisely, if $m\geq1$ and $\mathbf{H}\in \Gamma\bigl(T^{(k,m)}(TM)\otimes \mathbf{E}\bigr)$, then $\nabla \mathbf{H}\in \Gamma\bigl(T^{(k,m+1)}(TM)\otimes \mathbf{E}\bigr)$ and, in local coordinates, \[
\bigl(\operatorname{tr}_{\mathbf{g}}(\nabla \mathbf{H})\bigr)^{i_{1}\dots i_{k}a}_{j_{1}\dots j_{m-1}}
:=
g^{pq}\,(\nabla \mathbf{H})^{i_{1}\dots i_{k}a}_{pqj_{1}\dots j_{m-1}},
\] that is, the new derivative index is contracted with the first covariant index inherited from $\mathbf{H}$. \end{theorem}

Before proving the integration-by-parts formula, we recover several classical identities as special cases of Theorem~\ref{teo:green-m-esima-covariante-E}. These examples will clarify the meaning of each contraction and each boundary term.

\begin{example}[Classical integration-by-parts formulas as special cases of Theorem~\ref{teo:green-m-esima-covariante-E}]\label{ej:formulas-integracion-partes-casos-green-covariante} We show how several classical integration-by-parts formulas arise as special cases of Theorem~\ref{teo:green-m-esima-covariante-E}. In every part we take $s=1$ and $k=l=0$, so that \[
\mathbf{F}\in \Gamma_c(\mathbf{E}),
\qquad
\mathbf{G}\in \Gamma_c(T^{*}M\otimes \mathbf{E}).
\] When $\mathbf{E}$ is trivial of rank $1$, we identify $\Gamma_c(\mathbf{E})$ with $C_c^\infty(M)$ and $\Gamma_c(T^{*}M\otimes \mathbf{E})$ with $\Gamma_c(T^{*}M)$.

Recall that, if $(M,\mathbf{g})$ is a Riemannian manifold, the musical isomorphism $\flat\colon T_xM\longrightarrow T_x^*M$ defined by $X^\flat=\mathbf{g}(\mathbf{X},\cdot)$ is a linear isomorphism whose inverse is $\sharp\colon T_x^*M\longrightarrow T_xM$. Moreover, these isomorphisms preserve the inner product in the sense that, for $\alpha,\beta\in T_x^*M$ and $\mathbf{X},Y\in T_xM$, we have \[
\langle \alpha,\beta\rangle_{\mathbf{g}}
=
\mathbf{g}(\alpha^\sharp,\beta^\sharp),
\qquad
\mathbf{g}(\mathbf{X},Y)
=
\langle X^\flat,Y^\flat\rangle_{\mathbf{g}}.
\]

If $\mathbf{E}$ is a vector bundle equipped with a bundle metric $\mathbf{h}_{\mathbf{E}}$, recall that the inner product induced on $T^*M\otimes \mathbf{E}$ by $\mathbf{g}$ and $\mathbf{h}_{\mathbf{E}}$ is denoted by $\langle \cdot,\cdot\rangle_{\mathbf{g},\mathbf{h}_{\mathbf{E}}}$.

For simple tensor fields it is defined by \[
\langle \alpha\otimes s,\beta\otimes t\rangle_{\mathbf{g},\mathbf{h}_{\mathbf{E}}}
=
\langle \alpha,\beta\rangle_{\mathbf{g}}\, \mathbf{h}_{\mathbf{E}}(s,t).
\]

When $\mathbf{E}=M\times\mathbb{R}$ is trivial of rank $1$ with the usual metric $\mathbf{h}_{\mathbf{E}}(a,b)=ab$, there is a canonical identification \[
T^*M\otimes \mathbf{E} \cong T^*M,
\qquad
\alpha\otimes a \longleftrightarrow a\,\alpha.
\] Under this identification we have \[
\langle \alpha\otimes a,\beta\otimes b\rangle_{\mathbf{g},\mathbf{h}_{\mathbf{E}}}
=
ab\,\langle \alpha,\beta\rangle_{\mathbf{g}}
=
\langle a\alpha,b\beta\rangle_{\mathbf{g}},
\] so the induced product $\langle\cdot,\cdot\rangle_{\mathbf{g},\mathbf{h}_{\mathbf{E}}}$ agrees with the usual Riemannian inner product on $T^*M$. Consequently, when $\mathbf{E}$ is trivial of rank $1$, we omit the distinction between the two metrics.

\begin{enumerate}[label=\alph*)]

\item Let $\Omega\subset\mathbb{R}^{n}$ be a bounded open subset with smooth boundary. We consider $M=\overline{\Omega}$ with the Euclidean metric \[
\overline{\mathbf{g}}=\displaystyle\sum_{i=1}^{n}\mathbf{d}x^{i}\otimes \mathbf{d}x^{i},
\qquad
(\overline g_{ij})=\delta_{ij},\quad (\overline g^{ij})=\delta^{ij}.
\] We take the trivial bundle $\mathbf{E}=M\times\mathbb{R}$ with the usual bundle metric $\mathbf{h}_{\mathbf{E}}(a,b)=ab$ and the flat connection induced by the Euclidean Levi--Civita connection.

We apply Theorem~\ref{teo:green-m-esima-covariante-E}. The corresponding identity is \[
\int_{M}\langle \nabla \mathbf{F},\mathbf{G}\rangle_{\overline{\mathbf{g}},\mathbf{h}_{\mathbf{E}}}\,d\lambda_{n}
=
-\int_{M}\langle \mathbf{F},(\operatorname{tr}_{\overline{\mathbf{g}}}\circ \nabla)\mathbf{G}\rangle_{\mathbf{h}_{\mathbf{E}}}\,d\lambda_{n}
+
\int_{\partial M}\langle \mathbf{F},\iota_{\boldsymbol{\nu}}\mathbf{G}\rangle_{\mathbf{h}_{\mathbf{E}}}\,dS.
\]

Since $\mathbf{E}$ is trivial of rank $1$, this identity reduces to \[
\int_{M}\langle \nabla \mathbf{F},\mathbf{G}\rangle_{\overline{\mathbf{g}}}\,d\lambda_{n}
=
-\int_{M}\mathbf{F}\,(\operatorname{tr}_{\overline{\mathbf{g}}}\circ \nabla)\mathbf{G}\,d\lambda_{n}
+
\int_{\partial M}\mathbf{F}\,\iota_{\boldsymbol{\nu}}\mathbf{G}\,dS.
\]

We take $u,v\in C^\infty(M)$ and choose \[
\mathbf{F}=u,
\qquad
\mathbf{G}=\nabla v.
\]

In general, for a smooth function $v$ on a Riemannian manifold, the covariant derivative agrees with the differential: \[
\nabla v = dv = (\partial_i v)\,\mathbf{d}x^i.
\] Thus, in Cartesian coordinates, \[
(\nabla v)_i=\partial_i v.
\]

Moreover, \[
\langle \nabla u,\nabla v\rangle_{\overline{\mathbf{g}}}
=
\overline g^{ij}(\partial_i u)(\partial_j v)
=
\displaystyle\sum_{i=1}^n (\partial_i u)(\partial_i v).
\]

The general expression for the covariant derivative of a $1$-form is \[
(\nabla H)_{ij}
=
\partial_i H_j-\Gamma^k_{ij}H_k.
\] Applying this to $H=\nabla v$ gives \[
(\nabla(\nabla v))_{ij}
=
\partial_i\partial_j v-\Gamma^k_{ij}\partial_k v.
\] Since the Euclidean Levi--Civita connection satisfies $\Gamma^k_{ij}=0$, we obtain \[
(\operatorname{tr}_{\overline{\mathbf{g}}}\circ\nabla)(\nabla v)
=
\delta^{ij}\partial_i\partial_j v
=
\Delta v.
\]

Moreover, \[
\iota_{\boldsymbol{\nu}}(\nabla v)=\partial_{\boldsymbol{\nu}} v.
\]

Therefore, \[
\int_{\Omega}\nabla u\cdot\nabla v\,d\lambda_{n}
=
-\int_{\Omega}u\,\Delta v\,d\lambda_{n}
+
\int_{\partial\Omega}u\,\partial_{\boldsymbol{\nu}} v\,dS.
\]

\item Under the preceding hypotheses, subtracting the identities obtained for $(u,v)$ and $(v,u)$ gives \[
\int_{\Omega}\bigl(u\,\Delta v-v\,\Delta u\bigr)\,d\lambda_{n}
=
\int_{\partial\Omega}\bigl(u\,\partial_{\boldsymbol{\nu}} v-v\,\partial_{\boldsymbol{\nu}} u\bigr)\,dS,
\] which is known as Green's second identity.

\item Now let $(M,\mathbf{g})$ be a Riemannian manifold with boundary. Again take $\mathbf{E}=M\times\mathbb{R}$. For $u,v\in C_c^\infty(M)$, take $\mathbf{F}=u$ and $\mathbf{G}=\nabla v$, with $\nabla$ the Levi-Civita connection.

The abstract identity is \[
\int_M \langle \nabla u,\nabla v\rangle_{\mathbf{g},\mathbf{h}_{\mathbf{E}}}\,d\lambda_{\mathbf{g}}
=
-\int_M \langle u,(\operatorname{tr}_{\mathbf{g}}\circ\nabla)(\nabla v)\rangle_{\mathbf{h}_{\mathbf{E}}}\,d\lambda_{\mathbf{g}}
+
\int_{\partial M} \langle u,\iota_{\boldsymbol{\nu}}(\nabla v)\rangle_{\mathbf{h}_{\mathbf{E}}}\,d\lambda_{\widetilde{\mathbf{g}}}.
\]

Since $\mathbf{E}$ is trivial of rank $1$, it reduces to \[
\int_M \langle \nabla u,\nabla v\rangle_{\mathbf{g}}\,d\lambda_{\mathbf{g}}
=
-\int_M u\,(\operatorname{tr}_{\mathbf{g}}\circ\nabla)(\nabla v)\,d\lambda_{\mathbf{g}}
+
\int_{\partial M} u\,\iota_{\boldsymbol{\nu}}(\nabla v)\,d\lambda_{\widetilde{\mathbf{g}}}.
\]

In local coordinates, \[
(\nabla v)_j=\partial_j v,
\qquad
(\nabla(\nabla v))_{ij}
=
\partial_i\partial_j v-\Gamma^k_{ij}\partial_k v.
\] Therefore, \[
(\operatorname{tr}_{\mathbf{g}}\circ\nabla)(\nabla v)
=
g^{ij}\bigl(\partial_i\partial_j v-\Gamma^k_{ij}\partial_k v\bigr).
\] Let us verify that this expression agrees with the coordinate expression for the Laplace--Beltrami Laplacian, \[
\Delta_{\mathbf{g}} v
=
\frac{1}{\sqrt{\det(\mathbf{g})}}\,\partial_i\left(\sqrt{\det(\mathbf{g})}\,g^{ij}\partial_j v\right).
\]

By the definition of the metric trace, \[
(\operatorname{tr}_{\mathbf{g}}\circ\nabla)(\nabla v)
=
g^{ij}(\nabla(\nabla v))_{ij},
\] and since $(\nabla v)_j=\partial_j v$, we have \[
(\nabla(\nabla v))_{ij}
=
\partial_i\partial_j v-\Gamma^{k}_{ij}\partial_k v.
\] Therefore, \begin{equation}\label{eq:trace-hess}
(\operatorname{tr}_{\mathbf{g}}\circ\nabla)(\nabla v)
=
g^{ij}\partial_i\partial_j v
-
g^{ij}\Gamma^{k}_{ij}\partial_k v.
\end{equation}

On the other hand, \[
\Delta_{\mathbf{g}} v
=
\frac{1}{\sqrt{\det(\mathbf{g})}}\,\partial_i\!\bigl(\sqrt{\det(\mathbf{g})}\,g^{ij}\partial_j v\bigr).
\] Applying the Leibniz rule, \[
\Delta_{\mathbf{g}} v
=
\frac{1}{\sqrt{\det(\mathbf{g})}}
\Bigl(
\partial_i(\sqrt{\det(\mathbf{g})}\,g^{ij})\,\partial_j v
+
\sqrt{\det(\mathbf{g})}\,g^{ij}\partial_i\partial_j v
\Bigr)
\] \[
=
g^{ij}\partial_i\partial_j v
+
\frac{1}{\sqrt{\det(\mathbf{g})}}\partial_i(\sqrt{\det(\mathbf{g})}\,g^{ij})\,\partial_j v.
\]

We now expand the remaining term using the Leibniz rule once more: \[
\partial_i\left(\sqrt{\det(\mathbf{g})}\,g^{ij}\right)
=
\left(\partial_i\sqrt{\det(\mathbf{g})}\right)\,g^{ij}
+
\sqrt{\det(\mathbf{g})}\,\partial_i g^{ij}.
\] Dividing by $\sqrt{\det(\mathbf{g})}$ gives \[
\frac{1}{\sqrt{\det(\mathbf{g})}}\partial_i\left(\sqrt{\det(\mathbf{g})}\,g^{ij}\right)
=
g^{ij}\frac{\partial_i\sqrt{\det(\mathbf{g})}}{\sqrt{\det(\mathbf{g})}}
+
\partial_i g^{ij}.
\] By Remark~\ref{derivada raiz del determinante de la metrica}, \[
\frac{\partial_i\sqrt{\det(\mathbf{g})}}{\sqrt{\det(\mathbf{g})}}=\Gamma^{k}_{ik},
\] and, by compatibility of the connection with the metric (see \ref{derivada parcial metrica inversa}), \[
\partial_i g^{ij}+\Gamma^{k}_{ik}g^{ij}=-\,g^{ik}\Gamma^{j}_{ik},
\] we conclude that \[
\frac{1}{\sqrt{\det(\mathbf{g})}}\partial_i\left(\sqrt{\det(\mathbf{g})}\,g^{ij}\right)
=
-\,g^{ik}\Gamma^{j}_{ik}.
\] Substitution into the expression for $\Delta_{\mathbf{g}} v$ yields \[
\Delta_{\mathbf{g}} v
=
g^{ij}\partial_i\partial_j v
-
g^{ik}\Gamma^{j}_{ik}\partial_j v.
\] Renaming the dummy index $j\mapsto k$ gives \[
\Delta_{\mathbf{g}} v
=
g^{ij}\partial_i\partial_j v
-
g^{ij}\Gamma^{k}_{ij}\partial_k v,
\] which agrees with \eqref{eq:trace-hess}. Therefore, \[
(\operatorname{tr}_{\mathbf{g}}\circ\nabla)(\nabla v)=\Delta_{\mathbf{g}} v.
\]

Moreover, \[
\langle \nabla u,\nabla v\rangle_{\mathbf{g}}
=
g^{ij}(\partial_i u)(\partial_j v)
=
\mathbf{g}(\operatorname{grad}_{\mathbf{g}} u,\operatorname{grad}_{\mathbf{g}} v).
\]

Finally, \[
\iota_{\boldsymbol{\nu}}(\nabla v)
=
\nu^j\partial_j v
=
\mathbf{g}(\operatorname{grad}_{\mathbf{g}} v,\boldsymbol{\nu})
=
\partial_{\boldsymbol{\nu}} v.
\]

Substituting into the theorem gives \[
\int_M \mathbf{g}(\operatorname{grad}_{\mathbf{g}} u,\operatorname{grad}_{\mathbf{g}} v)\,d\lambda_{\mathbf{g}}
=
-\int_M u\,\Delta_{\mathbf{g}} v\,d\lambda _{\mathbf{g}}
+
\int_{\partial M} u\,\partial_{\boldsymbol{\nu}} v\,d\lambda_{\widetilde{\mathbf{g}}}.
\]

\item Let $\mathbf{X}\in \Gamma_c(TM)$ and let $\boldsymbol{\alpha}=\mathbf{X}^\flat$. We apply Theorem~\ref{teo:green-m-esima-covariante-E} with $\mathbf{F}\equiv 1$ and $\mathbf{G}=\boldsymbol{\alpha}$. Since $\nabla(1)=0$, we obtain \[
0
=
-\int_M (\operatorname{tr}_{\mathbf{g}}\circ\nabla)\boldsymbol{\alpha}\,d\lambda_{\mathbf{g}}
+
\int_{\partial M}\iota_{\boldsymbol{\nu}}\boldsymbol{\alpha}\,d\lambda_{\widetilde{\mathbf{g}}}.
\]

In local coordinates we have $(\nabla\boldsymbol{\alpha})_{ij}=\partial_i\alpha_j-\Gamma^k_{ij}\alpha_k$, and hence \[
(\operatorname{tr}_{\mathbf{g}}\circ\nabla)\boldsymbol{\alpha}
=
g^{ij}(\partial_i\alpha_j-\Gamma^k_{ij}\alpha_k).
\] Writing $\alpha_j=g_{jk}X^k$ and applying the Leibniz rule, \[
\partial_i\alpha_j
=
\partial_i(g_{jk}X^k)
=
(\partial_i g_{jk})X^k
+
g_{jk}\partial_i X^k.
\] Substituting into the preceding expression gives \[
(\operatorname{tr}_{\mathbf{g}}\circ\nabla)\boldsymbol{\alpha}
=
g^{ij}(\partial_i g_{jk})X^k
+
g^{ij}g_{jk}\partial_i X^k
-
g^{ij}\Gamma^m_{ij}g_{mk}X^k.
\] Since $g^{ij}g_{jk}=\delta^i_k$, the second term is $\partial_i X^i$.

To simplify the remaining terms, we use Corollary~\ref{compatibilidad levi civita metrica en coordenadas}, which in coordinates states that \[
\partial_i g_{jk}
=
\Gamma^s_{ij}g_{sk}
+
\Gamma^s_{ik}g_{js}.
\] Multiplying by $g^{ij}$ gives \[
g^{ij}\partial_i g_{jk}
=
g^{ij}\Gamma^s_{ij}g_{sk}
+
g^{ij}\Gamma^s_{ik}g_{js}.
\] The second term satisfies $g^{ij}\Gamma^s_{ik}g_{js}
=
\Gamma^s_{ik}\delta^i_s
=
\Gamma^i_{ik}$, and therefore \[
g^{ij}(\partial_i g_{jk})
-
g^{ij}\Gamma^s_{ij}g_{sk}
=
\Gamma^i_{ik}.
\] Substitution into the expression for $(\operatorname{tr}_{\mathbf{g}}\circ\nabla)\boldsymbol{\alpha}$ shows that \[
(\operatorname{tr}_{\mathbf{g}}\circ\nabla)\boldsymbol{\alpha}
=
\partial_i X^i+\Gamma^i_{ik}X^k
=
\operatorname{div}_{\mathbf{g}} \mathbf{X}.
\]

On the other hand, $\iota_{\boldsymbol{\nu}}\boldsymbol{\alpha}=\boldsymbol{\alpha}(\boldsymbol{\nu})=\mathbf{g}(\mathbf{X},\boldsymbol{\nu})$. Substituting into the initial identity finally gives \[
\int_M \operatorname{div}_{\mathbf{g}} \mathbf{X}\,d\lambda_{\mathbf{g}}
=
\int_{\partial M} \mathbf{g}(\mathbf{X},\boldsymbol{\nu})\,d\lambda_{\widetilde{\mathbf{g}}}.
\]

\end{enumerate}
\end{example}

The preceding examples illustrate the structure that the general formula must preserve. We now turn to its proof.

\begin{proof}[Proof of Theorem~\ref{teo:green-m-esima-covariante-E}] We proceed by induction on $s$. The general idea is to “transfer” derivatives from $\mathbf{F}$ to $\mathbf{G}$ by covariant integration by parts, keeping track of the boundary terms that arise at each step.

\medskip

\noindent\emph{Case $s=1$.} Fix a point $p\in\operatorname{Int}(M)$ and choose normal coordinates centered at $p$, with associated coordinate frame $(\boldsymbol{\partial}_{1},\dots,\boldsymbol{\partial}_{n})$. In these coordinates, \[
g_{ij}(p)=\delta_{ij},
\qquad
\Gamma^{k}_{ij}(p)=0,
\] and, in particular, $\nabla_{\boldsymbol{\partial}_{i}}\boldsymbol{\partial}_{j}(p)=0$ for all indices $i,j\in\{1,\dots,n\}$.

Now let $\mathbf{F}\in\Gamma_{c}\bigl(T^{(k,l)}(TM)\otimes \mathbf{E}\bigr)$ and $\mathbf{G}\in\Gamma_{c}\bigl(T^{(k,l+1)}(TM)\otimes \mathbf{E}\bigr)$. Define a map \[
\Phi\colon \mathfrak X(M)\longrightarrow C^{\infty}(M,\mathbb K),
\qquad
\Phi(\mathbf{Y}):=\big\langle \mathbf{F},\,\iota_{\mathbf{Y}}\mathbf{G}\big\rangle_{\mathbf{g},\mathbf{h}_{\mathbf{E}}}.
\] The map $\mathbf{Y}\mapsto \iota_{\mathbf{Y}}\mathbf{G}$ is $C^{\infty}(M,\mathbb R)$--linear in $\mathbf{Y}$, by the multilinearity of the tensor field $\mathbf{G}$ in its arguments. Consequently, $\Phi$ defines a smooth real $1$--form if \(\mathbf E\) is real and a smooth complex $1$--form if \(\mathbf E\) is complex. In the latter case, all the following arguments take place in the complexified bundles, and the divergence theorem is applied separately to the real and imaginary parts. Thus consider $\mathbf{X}:=\Phi^{\sharp}$, which defines a smooth vector field, complexified in the complex case, characterized (according to Definition~\ref{isomorfismos musicales}) by \[
\langle \mathbf{X},\mathbf{Y}\rangle_{\mathbf{g}}=\Phi(\mathbf{Y})=\big\langle \mathbf{F},\,\iota_{\mathbf{Y}}\mathbf{G}\big\rangle_{\mathbf{g},\mathbf{h}_{\mathbf{E}}},
\qquad
\forall\,\mathbf{Y}\in\mathfrak X(M).
\] In particular, for each $c\in\{1,\dots,n\}$, \begin{equation}\label{eq:X-parcial}
\langle \mathbf{X},\boldsymbol{\partial}_{c}\rangle_{\mathbf{g}}
=
\big\langle \mathbf{F},\,\iota_{\boldsymbol{\partial}_{c}}\mathbf{G}\big\rangle_{\mathbf{g},\mathbf{h}_{\mathbf{E}}}.
\end{equation}

At the point $p$, the divergence of $\mathbf{X}$ is expressed as \[
\operatorname{div}_{\mathbf{g}}\mathbf{X}(p)
=
\displaystyle\sum_{c=1}^{n}\big\langle \nabla_{\boldsymbol{\partial}_{c}}\mathbf{X},\boldsymbol{\partial}_{c}\big\rangle_{\mathbf{g}}(p).
\] Compatibility of $\nabla$ with $\mathbf{g}$ tells us that \[
\big\langle \nabla_{\boldsymbol{\partial}_{c}}\mathbf{X},\boldsymbol{\partial}_{c}\big\rangle_{\mathbf{g}}
=
\boldsymbol{\partial}_{c}\big(\langle \mathbf{X},\boldsymbol{\partial}_{c}\rangle_{\mathbf{g}}\big)
-
\big\langle \mathbf{X},\nabla_{\boldsymbol{\partial}_{c}}\boldsymbol{\partial}_{c}\big\rangle_{\mathbf{g}}.
\] Evaluating at $p$ and using $\nabla_{\boldsymbol{\partial}_{c}}\boldsymbol{\partial}_{c}(p)=0$ (since $\Gamma^{k}_{cc}(p)=0$), we find that \begin{equation}\label{eq:divX-en-p-derivada}
\operatorname{div}_{\mathbf{g}}\mathbf{X}(p)
=
\displaystyle\sum_{c=1}^{n}\boldsymbol{\partial}_{c}\big(\langle \mathbf{X},\boldsymbol{\partial}_{c}\rangle_{\mathbf{g}}\big)(p).
\end{equation} Substituting \eqref{eq:X-parcial} into \eqref{eq:divX-en-p-derivada} gives \begin{equation}\label{eq:divX-en-p-con-iota}
\operatorname{div}_{\mathbf{g}}\mathbf{X}(p)
=
\displaystyle\sum_{c=1}^{n}\boldsymbol{\partial}_{c}\Big(\big\langle \mathbf{F},\,\iota_{\boldsymbol{\partial}_{c}}\mathbf{G}\big\rangle_{\mathbf{g},\mathbf{h}_{\mathbf{E}}}\Big)(p).
\end{equation}

The induced connection on $T^{(k,l)}(TM)\otimes \mathbf{E}$ (which we continue to denote by $\nabla$) is compatible with the inner product $\langle\cdot,\cdot\rangle_{\mathbf{g},\mathbf{h}_{\mathbf{E}}}$. Consequently, for each $c$, \begin{equation}\label{eq:derivada-producto-interno}
\boldsymbol{\partial}_{c}\Big(\big\langle \mathbf{F},\,\iota_{\boldsymbol{\partial}_{c}}\mathbf{G}\big\rangle_{\mathbf{g},\mathbf{h}_{\mathbf{E}}}\Big)
=
\Big\langle \nabla_{\boldsymbol{\partial}_c}\mathbf{F},\,\iota_{\boldsymbol{\partial}_{c}}\mathbf{G}\Big\rangle_{\mathbf{g},\mathbf{h}_{\mathbf{E}}}
+
\Big\langle \mathbf{F},\,\nabla_{\boldsymbol{\partial}_c}\bigl(\iota_{\boldsymbol{\partial}_{c}}\mathbf{G}\bigr)\Big\rangle_{\mathbf{g},\mathbf{h}_{\mathbf{E}}}.
\end{equation} Evaluating \eqref{eq:derivada-producto-interno} at $p$ and substituting into \eqref{eq:divX-en-p-con-iota} gives \begin{equation}\label{eq:divX-en-p-dos-sumas-flat}
\operatorname{div}_{\mathbf{g}}\mathbf{X}(p)
=
\displaystyle\sum_{c=1}^{n}
\Big(
\big\langle \nabla_{\boldsymbol{\partial}_c}\mathbf{F},\,\iota_{\boldsymbol{\partial}_{c}}\mathbf{G}\big\rangle_{\mathbf{g},\mathbf{h}_{\mathbf{E}}}(p)
+
\big\langle \mathbf{F},\,\nabla_{\boldsymbol{\partial}_c}\bigl(\iota_{\boldsymbol{\partial}_{c}}\mathbf{G}\bigr)\big\rangle_{\mathbf{g},\mathbf{h}_{\mathbf{E}}}(p)
\Big).
\end{equation}

For the second term we use the identity (the Leibniz rule for contraction) \[
\nabla_{\boldsymbol{\partial}_c}\bigl(\iota_{\boldsymbol{\partial}_{c}}\mathbf{G}\bigr)
=
\iota_{\boldsymbol{\partial}_{c}}\bigl(\nabla_{\boldsymbol{\partial}_{c}}\mathbf{G}\bigr)
+
\iota_{\nabla_{\boldsymbol{\partial}_{c}}\boldsymbol{\partial}_{c}}\mathbf{G},
\] which follows from Lemma~\ref{leibniz multiplicacion interior}. Evaluating at $p$ and using $\nabla_{\boldsymbol{\partial}_{c}}\boldsymbol{\partial}_{c}(p)=0$ gives \[
\nabla_{\boldsymbol{\partial}_c}\bigl(\iota_{\boldsymbol{\partial}_{c}}\mathbf{G}\bigr)(p)
=
\iota_{\boldsymbol{\partial}_{c}}\bigl(\nabla_{\boldsymbol{\partial}_{c}}\mathbf{G}\bigr)(p).
\] Substituting this into \eqref{eq:divX-en-p-dos-sumas-flat} yields \begin{equation}\label{eq:divX-en-p-dos-sumas}
\operatorname{div}_{\mathbf{g}}\mathbf{X}(p)
=
\displaystyle\sum_{c=1}^{n}
\Big(
\big\langle \nabla_{\boldsymbol{\partial}_c}\mathbf{F},\,\iota_{\boldsymbol{\partial}_{c}}\mathbf{G}\big\rangle_{\mathbf{g},\mathbf{h}_{\mathbf{E}}}(p)
+
\big\langle \mathbf{F},\,\iota_{\boldsymbol{\partial}_{c}}(\nabla_{\boldsymbol{\partial}_c}\mathbf{G})\big\rangle_{\mathbf{g},\mathbf{h}_{\mathbf{E}}}(p)
\Big).
\end{equation}

In the chosen normal coordinates we have \[
 \nabla \mathbf{F}=\sum_{c=1}^n \mathbf{d}x^c\otimes\nabla_{\boldsymbol{\partial}_c}\mathbf{F},
 \qquad
 \iota_{\boldsymbol{\partial}_c}\mathbf{G}=\mathbf{G}(\boldsymbol{\partial}_c,\,\cdot\,),
\] because the derivative index and the inserted index both occupy the first position. The definition of the tensor product metric then implies \[
 \sum_{c=1}^n
 \big\langle\nabla_{\boldsymbol{\partial}_c}\mathbf{F},\iota_{\boldsymbol{\partial}_c}\mathbf{G}\big\rangle_{\mathbf{g},\mathbf{h}_{\mathbf{E}}}
 =\langle\nabla \mathbf{F},\mathbf{G}\rangle_{\mathbf{g},\mathbf{h}_{\mathbf{E}}}
 \quad\hbox{en }p.
\] Likewise, \[
 \operatorname{tr}_{\mathbf{g}}(\nabla \mathbf{G})
 =\sum_{c=1}^n(\nabla_{\boldsymbol{\partial}_c}\mathbf{G})(\boldsymbol{\partial}_c,\,\cdot\,)
 =\sum_{c=1}^n\iota_{\boldsymbol{\partial}_c}(\nabla_{\boldsymbol{\partial}_c}\mathbf{G})
 \quad\hbox{en }p,
\] and the second sum in \eqref{eq:divX-en-p-dos-sumas} is $\langle \mathbf{F},\operatorname{tr}_{\mathbf{g}}(\nabla \mathbf{G})\rangle_{\mathbf{g},\mathbf{h}_{\mathbf{E}}}(p)$.

\medskip

Combining the two identifications gives \[
\operatorname{div}_{\mathbf{g}}\mathbf{X}(p)
=
\big\langle \nabla \mathbf{F},\mathbf{G}\big\rangle_{\mathbf{g},\mathbf{h}_{\mathbf{E}}}(p)
+
\big\langle \mathbf{F},\operatorname{tr}_{\mathbf{g}}(\nabla \mathbf{G})\big\rangle_{\mathbf{g},\mathbf{h}_{\mathbf{E}}}(p).
\] Since $p\in\operatorname{Int}(M)$ is arbitrary, the pointwise identity holds on \(\operatorname{Int}(M)\). Both sides are smooth functions up to the boundary, and \(\operatorname{Int}(M)\) is dense in \(M\); by continuity, the identity also holds on \(\partial M\). It therefore holds throughout $M$: \[
\operatorname{div}_{\mathbf{g}}\mathbf{X}
=
\langle \nabla \mathbf{F},\mathbf{G}\rangle_{\mathbf{g},\mathbf{h}_{\mathbf{E}}}
+
\langle \mathbf{F},\operatorname{tr}_{\mathbf{g}}(\nabla \mathbf{G})\rangle_{\mathbf{g},\mathbf{h}_{\mathbf{E}}}.
\]

Integrating over $M$ and applying Theorem~\ref{teo:divergencia} (the divergence theorem) gives \[
\int_{M}\operatorname{div}_{\mathbf{g}}\mathbf{X}\,d\lambda_{\mathbf{g}}
=
\int_{\partial M}\langle \mathbf{X},\boldsymbol{\nu}\rangle_{\mathbf{g}}\,d\lambda_{\widetilde{\mathbf{g}}}.
\] The field \(\mathbf X\) has compact support, since it vanishes wherever \(\mathbf F\) or \(\mathbf G\) vanishes; thus the divergence theorem also applies when \(M\) is noncompact. If \(\partial M=\varnothing\), the right-hand side is zero by the convention in the statement. On the boundary, by the definition of $\mathbf{X}$ (and of $\Phi$), \[
\langle \mathbf{X},\boldsymbol{\nu}\rangle_{\mathbf{g}}
=
(\mathbf{X}^{\flat})(\boldsymbol{\nu})
=
\Phi(\boldsymbol{\nu})
=
\big\langle \mathbf{F},\,\iota_{\boldsymbol{\nu}}\mathbf{G}\big\rangle_{\mathbf{g},\mathbf{h}_{\mathbf{E}}}.
\] Combining this with the preceding pointwise identity gives \[
\int_{M}\big\langle \nabla \mathbf{F},\mathbf{G}\big\rangle_{\mathbf{g},\mathbf{h}_{\mathbf{E}}}\,d\lambda_{\mathbf{g}}
=
-\int_{M}\big\langle \mathbf{F},\operatorname{tr}_{\mathbf{g}}(\nabla \mathbf{G})\big\rangle_{\mathbf{g},\mathbf{h}_{\mathbf{E}}}\,d\lambda_{\mathbf{g}}
+
\int_{\partial M}\big\langle \mathbf{F},\,\iota_{\boldsymbol{\nu}}\mathbf{G}\big\rangle_{\mathbf{g},\mathbf{h}_{\mathbf{E}}}\,d\lambda_{\widetilde{\mathbf{g}}}.
\] This proves the case $s=1$.

\medskip

We now verify the induction step. Suppose that the identity in the statement holds for some integer $s\ge 1$, that is, for every \[
\mathbf{F}\in\Gamma_{c}\bigl(T^{(k,l)}(TM)\otimes \mathbf{E}\bigr),
\quad
\mathbf{H}\in\Gamma_{c}\bigl(T^{(k,l+s)}(TM)\otimes \mathbf{E}\bigr),
\] we have \[
\int_{M}\langle \nabla^{s}\mathbf{F},\mathbf{H}\rangle_{\mathbf{g},\mathbf{h}_{\mathbf{E}}}\,d\lambda_{\mathbf{g}}
=
(-1)^{s}\int_{M}\big\langle \mathbf{F},(\operatorname{tr}_{\mathbf{g}}\circ\nabla)^{s}\mathbf{H}\big\rangle_{\mathbf{g},\mathbf{h}_{\mathbf{E}}}\,d\lambda_{\mathbf{g}}
\] \[
\qquad
+
\displaystyle\sum_{j=0}^{s-1}(-1)^{\,s-1-j}
\int_{\partial M}\Big\langle \nabla^{j}\mathbf{F},\,\iota_{\boldsymbol{\nu}}\bigl((\operatorname{tr}_{\mathbf{g}}\circ\nabla)^{\,s-1-j}\mathbf{H}\bigr)\Big\rangle_{\mathbf{g},\mathbf{h}_{\mathbf{E}}}\,d\lambda_{\widetilde{\mathbf{g}}}.
\]

Our goal is to show that the same identity holds for $s+1$. Thus let \[
\mathbf{G}\in\Gamma_{c}\bigl(T^{(k,l+s+1)}(TM)\otimes \mathbf{E}\bigr).
\] We apply the case $s=1$ to the pair of fields $\nabla^{s}\mathbf{F}$ and $\mathbf{G}$. This gives \[
\int_{M}\langle \nabla^{s+1}\mathbf{F},\mathbf{G}\rangle_{\mathbf{g},\mathbf{h}_{\mathbf{E}}}\,d\lambda_{\mathbf{g}}
=
-\int_{M}\big\langle \nabla^{s}\mathbf{F},\operatorname{tr}_{\mathbf{g}}(\nabla \mathbf{G})\big\rangle_{\mathbf{g},\mathbf{h}_{\mathbf{E}}}\,d\lambda_{\mathbf{g}}
+
\int_{\partial M}\big\langle \nabla^{s}\mathbf{F},\,\iota_{\boldsymbol{\nu}}\mathbf{G}\big\rangle_{\mathbf{g},\mathbf{h}_{\mathbf{E}}}\,d\lambda_{\widetilde{\mathbf{g}}}.
\]

Now define \[
\mathbf{H}:=\operatorname{tr}_{\mathbf{g}}(\nabla \mathbf{G})
\in\Gamma_{c}\bigl(T^{(k,l+s)}(TM)\otimes \mathbf{E}\bigr).
\] Applying the induction hypothesis to the pair $(\mathbf{F},\mathbf{H})$ gives \[
\int_{M}\langle \nabla^{s}\mathbf{F},\mathbf{H}\rangle_{\mathbf{g},\mathbf{h}_{\mathbf{E}}}\,d\lambda_{\mathbf{g}}
=
(-1)^{s}
\int_{M}\big\langle \mathbf{F},(\operatorname{tr}_{\mathbf{g}}\circ\nabla)^{s}\mathbf{H}\big\rangle_{\mathbf{g},\mathbf{h}_{\mathbf{E}}}\,d\lambda_{\mathbf{g}}
\] \[
\qquad
+
\displaystyle\sum_{j=0}^{s-1}(-1)^{\,s-1-j}
\int_{\partial M}\Big\langle \nabla^{j}\mathbf{F},\,\iota_{\boldsymbol{\nu}}\bigl((\operatorname{tr}_{\mathbf{g}}\circ\nabla)^{\,s-1-j}\mathbf{H}\bigr)\Big\rangle_{\mathbf{g},\mathbf{h}_{\mathbf{E}}}\,d\lambda_{\widetilde{\mathbf{g}}}.
\]

Since $\mathbf{H}=\operatorname{tr}_{\mathbf{g}}(\nabla \mathbf{G})$, we have \[
(\operatorname{tr}_{\mathbf{g}}\circ\nabla)^{s}\mathbf{H}
=
(\operatorname{tr}_{\mathbf{g}}\circ\nabla)^{s+1}\mathbf{G},
\qquad
(\operatorname{tr}_{\mathbf{g}}\circ\nabla)^{\,s-1-j}\mathbf{H}
=
(\operatorname{tr}_{\mathbf{g}}\circ\nabla)^{\,s-j}\mathbf{G},
\] and therefore \[
\int_{M}\big\langle \nabla^{s}\mathbf{F},\operatorname{tr}_{\mathbf{g}}(\nabla \mathbf{G})\big\rangle_{\mathbf{g},\mathbf{h}_{\mathbf{E}}}\,d\lambda_{\mathbf{g}}
=
(-1)^{s}
\int_{M}\big\langle \mathbf{F},(\operatorname{tr}_{\mathbf{g}}\circ\nabla)^{s+1}\mathbf{G}\big\rangle_{\mathbf{g},\mathbf{h}_{\mathbf{E}}}\,d\lambda_{\mathbf{g}}
\] \[
\qquad
+
\displaystyle\sum_{j=0}^{s-1}(-1)^{\,s-1-j}
\int_{\partial M}\Big\langle \nabla^{j}\mathbf{F},\,\iota_{\boldsymbol{\nu}}\bigl((\operatorname{tr}_{\mathbf{g}}\circ\nabla)^{\,s-j}\mathbf{G}\bigr)\Big\rangle_{\mathbf{g},\mathbf{h}_{\mathbf{E}}}\,d\lambda_{\widetilde{\mathbf{g}}}.
\]

Substituting this identity into the expression for \(\displaystyle \int_{M}\langle \nabla^{s+1}\mathbf{F},\mathbf{G}\rangle_{\mathbf{g},\mathbf{h}_{\mathbf{E}}}\,d\lambda_{\mathbf{g}}\) gives \[
\int_{M}\langle \nabla^{s+1}\mathbf{F},\mathbf{G}\rangle_{\mathbf{g},\mathbf{h}_{\mathbf{E}}}\,d\lambda_{\mathbf{g}}
=
(-1)^{s+1}
\int_{M}\big\langle \mathbf{F},(\operatorname{tr}_{\mathbf{g}}\circ\nabla)^{s+1}\mathbf{G}\big\rangle_{\mathbf{g},\mathbf{h}_{\mathbf{E}}}\,d\lambda_{\mathbf{g}}
\] \[
\qquad
+
\displaystyle\sum_{j=0}^{s-1}(-1)^{\,s-j}
\int_{\partial M}\Big\langle \nabla^{j}\mathbf{F},\,\iota_{\boldsymbol{\nu}}\bigl((\operatorname{tr}_{\mathbf{g}}\circ\nabla)^{\,s-j}\mathbf{G}\bigr)\Big\rangle_{\mathbf{g},\mathbf{h}_{\mathbf{E}}}\,d\lambda_{\widetilde{\mathbf{g}}}
+
\int_{\partial M}\big\langle \nabla^{s}\mathbf{F},\,\iota_{\boldsymbol{\nu}}\mathbf{G}\big\rangle_{\mathbf{g},\mathbf{h}_{\mathbf{E}}}\,d\lambda_{\widetilde{\mathbf{g}}}.
\]

The last boundary term already has the desired form if we interpret $(\operatorname{tr}_{\mathbf{g}}\circ\nabla)^0$ as the identity: \[
\int_{\partial M}\big\langle \nabla^{s}\mathbf{F},\,\iota_{\boldsymbol{\nu}}\mathbf{G}\big\rangle_{\mathbf{g},\mathbf{h}_{\mathbf{E}}}\,d\lambda_{\widetilde{\mathbf{g}}}
=
(-1)^{\,s-s}
\int_{\partial M}\Big\langle \nabla^{s}\mathbf{F},\,\iota_{\boldsymbol{\nu}}\bigl((\operatorname{tr}_{\mathbf{g}}\circ\nabla)^{\,s-s}\mathbf{G}\bigr)\Big\rangle_{\mathbf{g},\mathbf{h}_{\mathbf{E}}}\,d\lambda_{\widetilde{\mathbf{g}}}.
\] We can therefore incorporate this term into the preceding sum by taking $j=s$, obtaining \[
\int_{M}\langle \nabla^{s+1}\mathbf{F},\mathbf{G}\rangle_{\mathbf{g},\mathbf{h}_{\mathbf{E}}}\,d\lambda_{\mathbf{g}}
=
(-1)^{s+1}
\int_{M}\big\langle \mathbf{F},(\operatorname{tr}_{\mathbf{g}}\circ\nabla)^{s+1}\mathbf{G}\big\rangle_{\mathbf{g},\mathbf{h}_{\mathbf{E}}}\,d\lambda_{\mathbf{g}}
\] \[
\qquad
+
\displaystyle\sum_{j=0}^{s}(-1)^{\,s-j}
\int_{\partial M}\Big\langle \nabla^{j}\mathbf{F},\,\iota_{\boldsymbol{\nu}}\bigl((\operatorname{tr}_{\mathbf{g}}\circ\nabla)^{\,s-j}\mathbf{G}\bigr)\Big\rangle_{\mathbf{g},\mathbf{h}_{\mathbf{E}}}\,d\lambda_{\widetilde{\mathbf{g}}},
\] which is precisely the formula in the statement for $s+1$. The result then follows by induction on $s$. \end{proof}

As a corollary, we can extend the usual Green's formula to the Bochner Laplacian:

\begin{corollary}[Green's formula for the Bochner Laplacian]\label{cor:green-bochner}\index{Green's formula for the Bochner Laplacian@Green's formula for the Bochner Laplacian} Let $(M,\mathbf{g})$ be a Riemannian manifold with or without boundary, and let $\pi_{\mathbf{E}}\colon \mathbf{E}\longrightarrow M$ be a smooth vector bundle equipped with a bundle metric $\mathbf{h}_{\mathbf{E}}$ and a compatible connection $\nabla^{\mathbf{E}}$. Let $\boldsymbol{\nu}$ be the outward unit normal to $\partial M$.

Then, for any $\mathbf{u},\mathbf{v}\in \Gamma_{c}(\mathbf{E})$, we have \[
\int_{M}\langle \nabla^{\mathbf{E}}\mathbf{u},\nabla^{\mathbf{E}}\mathbf{v}\rangle_{\mathbf{g},\mathbf{h}_{\mathbf{E}}}\,d\lambda_{\mathbf{g}}
=
\int_{M}\langle \mathbf{u},\Delta_{B}\mathbf{v}\rangle_{\mathbf{h}_{\mathbf{E}}}\,d\lambda_{\mathbf{g}}
+
\int_{\partial M}\big\langle \mathbf{u},\nabla^{\mathbf{E}}_{\boldsymbol{\nu}}\mathbf{v}\big\rangle_{\mathbf{h}_{\mathbf{E}}}\,d\lambda_{\widetilde{\mathbf{g}}},
\] where, by definition, \[
\Delta_{B}:=(\nabla^{\mathbf{E}})_h^*\nabla^{\mathbf{E}}.
\] When $\partial M\neq\varnothing$, we use the formal differential adjoint up to the boundary from Remark~\ref{obs:adjunto-formal-frontera-dominio}; no domain for $L^2(M,\mathbf E)$ with boundary conditions has yet been chosen. \end{corollary}

\begin{proof} We apply Theorem~\ref{teo:green-m-esima-covariante-E} with $s=1$ and $k=l=0$, taking \[
\mathbf{F}=\mathbf{u}\in\Gamma_c(\mathbf{E}),
\qquad
\mathbf{G}=\nabla^{\mathbf{E}}\mathbf{v}\in \Gamma_c(T^{*}M\otimes \mathbf{E}).
\] The identity in the theorem takes the form \[
\int_{M}\langle \nabla^{\mathbf{E}}\mathbf{u},\nabla^{\mathbf{E}}\mathbf{v}\rangle_{\mathbf{g},\mathbf{h}_{\mathbf{E}}}\,d\lambda_{\mathbf{g}}
=
-\int_{M}\big\langle \mathbf{u},(\operatorname{tr}_{\mathbf{g}}\circ \nabla)(\nabla^{\mathbf{E}}\mathbf{v})\big\rangle_{\mathbf{h}_{\mathbf{E}}}\,d\lambda_{\mathbf{g}}
+
\int_{\partial M}\Big\langle \mathbf{u},\,\iota_{\boldsymbol{\nu}}(\nabla^{\mathbf{E}}\mathbf{v})\Big\rangle_{\mathbf{h}_{\mathbf{E}}}\,d\lambda_{\widetilde{\mathbf{g}}}.
\] By the convention adopted in the definition of the Bochner Laplacian, \[
\Delta_{B}=-\,(\operatorname{tr}_{\mathbf{g}}\circ \nabla)\circ \nabla^{\mathbf{E}},
\] and the interior term can therefore be rewritten as \[
-\big\langle \mathbf{u},(\operatorname{tr}_{\mathbf{g}}\circ \nabla)(\nabla^{\mathbf{E}}\mathbf{v})\big\rangle_{\mathbf{h}_{\mathbf{E}}}
=
\langle \mathbf{u},\Delta_{B}\mathbf{v}\rangle_{\mathbf{h}_{\mathbf{E}}}.
\]

It remains only to identify the boundary term. In local coordinates $(x^{1},\dots,x^{n})$ we write $\boldsymbol{\nu}=\nu^{i}\boldsymbol{\partial}_i$. By the definition of interior multiplication $\iota_{\boldsymbol{\nu}}$ in the first covariant index, $\iota_{\boldsymbol{\nu}}(\nabla^{\mathbf{E}}\mathbf{v})
=
\nu^{i}(\nabla^{\mathbf{E}}\mathbf{v})_i$. Moreover, $(\nabla^{\mathbf{E}}\mathbf{v})_i=\nabla^{\mathbf{E}}_{\boldsymbol{\partial}_i}\mathbf{v}$. Substitution gives $\iota_{\boldsymbol{\nu}}(\nabla^{\mathbf{E}}\mathbf{v})
=
\nu^{i}\nabla^{\mathbf{E}}_{\boldsymbol{\partial}_i}\mathbf{v}$. Finally, using the $C^{\infty}(M)$-linearity of the connection in the vector field, \[
\nabla^{\mathbf{E}}_{\boldsymbol{\nu}}\mathbf{v}
=
\nabla^{\mathbf{E}}_{\nu^{i}\boldsymbol{\partial}_i}\mathbf{v}
=
\nu^{i}\nabla^{\mathbf{E}}_{\boldsymbol{\partial}_i}\mathbf{v},
\] we conclude that $\iota_{\boldsymbol{\nu}}(\nabla^{\mathbf{E}}\mathbf{v})=\nabla^{\mathbf{E}}_{\boldsymbol{\nu}}\mathbf{v}$. Substituting into the preceding formula gives the stated identity. \end{proof}

The following result is another corollary of this integration-by-parts formula:

\begin{corollary}[Local expression for the formal adjoint of $\nabla^{s}$]\label{cor: expresion coordenada adjunto formal s derivada covariante}\index{local expression for the formal adjoint@local expression for the formal adjoint} Let $(M,\mathbf{g})$ be a Riemannian manifold without boundary, and let $\mathbf{E}\to M$ be a vector bundle with a bundle metric $\mathbf{h}_{\mathbf{E}}$ and a compatible connection $\nabla^{\mathbf{E}}$. For $s\ge 1$, the formal adjoint of \[
\nabla^{s}\colon \Gamma(\mathbf{E})\longrightarrow \Gamma\bigl(T^{(0,s)}(TM)\otimes \mathbf{E}\bigr)
\] is given by \[
(\nabla^{s})^{*}=(-1)^{s}(\operatorname{tr}_{\mathbf{g}}\circ\nabla)^{s}.
\] Moreover, in a chart $(U,x^{1},\dots,x^{n})$, if $\mathbf{G}\in \Gamma\bigl(T^{(0,s)}(TM)\otimes \mathbf{E}\bigr)$ has components $G^{a}_{j_{1}\dots j_{s}}$, then \begin{equation}\label{eq:adjunto formal s derivada covariante coordenadas}
\bigl((\nabla^{s})^{*}\mathbf{G}\bigr)^{a}
=
(-1)^{s}\Bigg[
\frac{1}{\sqrt{\det(\mathbf{g})}}\,
\partial_{i_{1}}\cdots\partial_{i_{s}}
\Big(\sqrt{\det(\mathbf{g})}\,g^{i_{1}j_{1}}\cdots g^{i_{s}j_{s}}\,G^{a}_{j_{1}\dots j_{s}}\Big)
+
\displaystyle\sum_{b=1}^{r}\displaystyle\sum_{|\beta|\le s-1}
\bigl(C_{\beta}\bigr)^{a\,j_{1}\dots j_{s}}_{b}\,\partial^{\beta}\!\bigl(G^{b}_{j_{1}\dots j_{s}}\bigr)
\Bigg],
\end{equation} where the coefficients $\bigl(C_{\beta}\bigr)^{a\,j_{1}\dots j_{s}}_{b}$ are smooth functions on $U$ that depend only on $g_{ij}$, the connection coefficients $\Gamma^{k}_{ij}$ and $(A_{i})^{a}_{b}$, and their partial derivatives of orders up to $s-1$. \end{corollary}

\begin{proof} Since $M$ has no boundary, the boundary terms in Theorem~\ref{teo:green-m-esima-covariante-E} vanish, and therefore \[
\int_{M}\langle \nabla^{s}\mathbf{F},\mathbf{G}\rangle_{\mathbf{g},\mathbf{h}_{\mathbf{E}}}\,d\lambda_{\mathbf{g}}
=
(-1)^{s}\int_{M}\big\langle \mathbf{F},(\operatorname{tr}_{\mathbf{g}}\circ \nabla)^{s}\mathbf{G}\big\rangle_{\mathbf{g},\mathbf{h}_{\mathbf{E}}}\,d\lambda_{\mathbf{g}},
\qquad
\forall\,\mathbf{F}\in\Gamma_{c}(\mathbf{E}),\ \forall\,\mathbf{G}\in\Gamma_c\bigl(T^{(0,s)}(TM)\otimes \mathbf{E}\bigr).
\] By Theorem~\ref{cor: adjunto global en coordenadas locales (existencia y expresion)} (existence and uniqueness of the formal adjoint), \[
(\nabla^{s})^{*}=(-1)^{s}(\operatorname{tr}_{\mathbf{g}}\circ\nabla)^{s}.
\] It remains to verify the local formula \eqref{eq:adjunto formal s derivada covariante coordenadas}. Fix a chart $(U,x^{1},\dots,x^{n})$ and a local frame $(\mathbf{e}_{1},\dots,\mathbf{e}_{r})$ of $\mathbf{E}\restriction_{U}$. Let $m\ge 0$ and $\mathbf{H}\in \Gamma\bigl(T^{(0,m+1)}(TM)\otimes \mathbf{E}\bigr)$ with components $H^{a}_{qj_{1}\dots j_{m}}$. By definition, \[
((\operatorname{tr}_{\mathbf{g}}\circ\nabla) \mathbf{H})^{a}_{j_{1}\dots j_{m}}
=
g^{pq}\,(\nabla_{\boldsymbol{\partial}_{p}}\mathbf{H})^{a}_{qj_{1}\dots j_{m}}.
\] Applying Lemma~\ref{lema: derivada covariante en coordenadas haz tensorial y E} (with $k=0$, $l=m+1$) gives \begin{equation}\label{eq:D-raw}
\begin{aligned}
((\operatorname{tr}_{\mathbf{g}}\circ\nabla)\mathbf{H})^{a}_{j_{1}\dots j_{m}}
&=
g^{pq}\partial_{p}\bigl(H^{a}_{qj_{1}\dots j_{m}}\bigr)
-\Gamma^{\alpha}_{p q}\,g^{pq}\,H^{a}_{\alpha j_{1}\dots j_{m}}\\
&\quad
-\displaystyle\sum_{\ell=1}^{m}\Gamma^{\alpha}_{p j_{\ell}}\,g^{pq}\,
H^{a}_{qj_{1}\dots j_{\ell-1}\alpha j_{\ell+1}\dots j_{m}}
+\displaystyle\sum_{b=1}^{r}(A_{p})^{a}_{b}\,g^{pq}H^{b}_{qj_{1}\dots j_{m}}.
\end{aligned}
\end{equation}

Using Remark~\ref{derivada raiz del determinante de la metrica}, which tells us that $\Gamma_{\alpha p}^{\alpha}=\displaystyle\frac{1}{\sqrt{\det(\mathbf{g})}}\partial_{p}\left(\sqrt{\det(\mathbf{g})}\right)$ and the identity in Corollary~\ref{derivada parcial metrica inversa}, which can be written as $\partial_{p}(g^{pq})=-\Gamma^{p}_{p\alpha}g^{\alpha q}-\Gamma^{q}_{p\alpha}g^{p\alpha},$ we obtain the equality \begin{equation}\label{eq:key-density-identity}
g^{pq}\partial_{p}\bigl(Y_{q}\bigr)-\Gamma^{\alpha}_{pq}g^{pq}Y_{\alpha}
=
\frac{1}{\sqrt{\det(\mathbf{g})}}\,\partial_{p}\Big(\sqrt{\det(\mathbf{g})}\,g^{pq}Y_{q}\Big),
\qquad
\forall (Y_{q})_{q=1}^{n}.
\end{equation} To verify it, we expand the right-hand side: \[
\frac{1}{\sqrt{\det(\mathbf{g})}}\partial_{p}\left(\sqrt{\det(\mathbf{g})}g^{pq} Y_{q}\right) = \frac{1}{\sqrt{\det(\mathbf{g})}}\left(\partial_{p}\left(\sqrt{\det(\mathbf{g})}\right)g^{pq}Y_{q}+\cancel{\sqrt{\det(\mathbf{g})}}\partial_{p}(g^{pq})Y_{q}+\cancel{\sqrt{\det(\mathbf{g})}}g^{pq}\partial_{p}(Y_{q})\right)
\] \[
= \Gamma_{\alpha p}^{\alpha}g^{pq}Y_{q}+{\color{red}\partial_{p}(g^{pq})Y_{q}}+g^{pq}\partial_{p}(Y_{q})
\] \[
=\underbrace{\Gamma_{p\alpha}^{p}g^{\alpha q}Y_{q}}_{p\leftrightarrow \alpha}{\color{red}-\Gamma_{p\alpha}^{p}g^{\alpha q}Y_{q}-\Gamma_{p\alpha }^{q}g^{p\alpha}Y_{q}}+g^{pq}\partial_{p}(Y_{q})
\] \[
= g^{pq}\partial_{p}(Y_{q})-\underbrace{\Gamma_{p\alpha}^{q}g^{p\alpha}Y_{q}}_{\alpha\leftrightarrow q}
\] \[
=g^{pq}\partial_{p}(Y_{q})-\Gamma_{pq}^{\alpha}g^{pq}Y_{\alpha}.
\]

Applying \eqref{eq:key-density-identity} with $Y_{q}=H^{a}_{qj_{1}\dots j_{m}}$, we can rewrite equation \eqref{eq:D-raw} as \begin{equation}\label{eq:D-local}
\begin{aligned}
((\operatorname{tr}_{\mathbf{g}}\circ\nabla)\mathbf{H})^{a}_{j_{1}\dots j_{m}}
&=\frac{1}{\sqrt{\det(\mathbf{g})}}\,\partial_{p}\Big(\sqrt{\det(\mathbf{g})}\,g^{pq}H^{a}_{qj_{1}\dots j_{m}}\Big)\\
&\quad-\displaystyle\sum_{\ell=1}^{m}\Gamma^{\alpha}_{p j_{\ell}}\,g^{pq}\,
H^{a}_{qj_{1}\dots j_{\ell-1}\alpha j_{\ell+1}\dots j_{m}}\\
&\quad+\displaystyle\sum_{b=1}^{r}(A_{p})^{a}_{b}\,g^{pq}H^{b}_{qj_{1}\dots j_{m}}.
\end{aligned}
\end{equation}

We prove by induction on $s\ge 1$ that, for every $\mathbf{G}\in\Gamma\bigl(T^{(0,s)}(TM)\otimes \mathbf{E}\bigr)$ with components $G^{a}_{j_{1}\dots j_{s}}$, we have \begin{equation}\label{eq:claim}
\bigl((\operatorname{tr}_{\mathbf{g}}\circ\nabla)^{s}\mathbf{G}\bigr)^{a}
=
\frac{1}{\sqrt{\det(\mathbf{g})}}\,
\partial_{i_{1}}\cdots\partial_{i_{s}}
\Big(\sqrt{\det(\mathbf{g})}\,g^{i_{1}j_{1}}\cdots g^{i_{s}j_{s}}\,G^{a}_{j_{1}\dots j_{s}}\Big)
+
\displaystyle\sum_{b=1}^{r}\displaystyle\sum_{|\beta|\le s-1}
\bigl(C_{\beta}\bigr)^{a\,j_{1}\dots j_{s}}_{b}\,\partial^{\beta}\!\bigl(G^{b}_{j_{1}\dots j_{s}}\bigr).
\end{equation}

For $s=1$, taking $m=0$ in \eqref{eq:D-local} gives \[
((\operatorname{tr}_{\mathbf{g}}\circ\nabla)\mathbf{G})^{a}
=
\frac{1}{\sqrt{\det(\mathbf{g})}}\,\partial_{p}\Big(\sqrt{\det(\mathbf{g})}\,g^{pq}G^{a}_{q}\Big)
+\displaystyle\sum_{b=1}^{r}(A_{p})^{a}_{b}\,g^{pq}G^{b}_{q},
\] which has the form \eqref{eq:claim} with $|\beta|= 0$. Now suppose that \eqref{eq:claim} holds for some $s\ge 1$, and let us prove the case $s+1$. Let $\mathbf{G}\in\Gamma\bigl(T^{(0,s+1)}(TM)\otimes \mathbf{E}\bigr)$ and define \[
\mathbf{H}:=(\operatorname{tr}_{\mathbf{g}}\circ\nabla)\mathbf{G}\in\Gamma\bigl(T^{(0,s)}(TM)\otimes \mathbf{E}\bigr).
\] By the induction hypothesis applied to $\mathbf{H}$, we have \begin{equation}\label{eq:IH-on-H}
\bigl((\operatorname{tr}_{\mathbf{g}}\circ\nabla)^{s}\mathbf{H}\bigr)^{a}
=
\frac{1}{\sqrt{\det(\mathbf{g})}}\,
\partial_{i_{1}}\cdots\partial_{i_{s}}
\Big(\sqrt{\det(\mathbf{g})}\,g^{i_{1}j_{1}}\cdots g^{i_{s}j_{s}}\,H^{a}_{j_{1}\dots j_{s}}\Big)
+
\displaystyle\sum_{b=1}^{r}\displaystyle\sum_{|\beta|\le s-1}
\bigl(C_{\beta}\bigr)^{a\,j_{1}\dots j_{s}}_{b}\,\partial^{\beta}\!\bigl(H^{b}_{j_{1}\dots j_{s}}\bigr).
\end{equation} Since, by definition, \[
(\operatorname{tr}_{\mathbf{g}}\circ\nabla)^{s}\mathbf{H}
=
(\operatorname{tr}_{\mathbf{g}}\circ\nabla)^{s+1}\mathbf{G},
\] it suffices to expand the components $H^{a}_{j_{1}\dots j_{s}}$ explicitly in terms of $\mathbf{G}$. Using \eqref{eq:D-local} with $m=s$ gives \begin{equation}\label{eq:H-expanded}
\begin{aligned}
H^{a}_{j_{1}\dots j_{s}}
&=
\frac{1}{\sqrt{\det(\mathbf{g})}}\,\partial_{p}\Big(\sqrt{\det(\mathbf{g})}\,g^{pq}G^{a}_{qj_{1}\dots j_{s}}\Big)
-\displaystyle\sum_{\ell=1}^{s}\Gamma^{\alpha}_{p j_{\ell}}\,g^{pq}\,
G^{a}_{qj_{1}\dots j_{\ell-1}\alpha j_{\ell+1}\dots j_{s}}\\
&\quad
+\displaystyle\sum_{b=1}^{r}(A_{p})^{a}_{b}\,g^{pq}G^{b}_{qj_{1}\dots j_{s}}.
\end{aligned}
\end{equation}

Substituting \eqref{eq:H-expanded} into the leading term of \eqref{eq:IH-on-H} and separating the summands gives \begin{equation}\label{eq:split-I-II-III}
\frac{1}{\sqrt{\det(\mathbf{g})}}\,
\partial_{i_{1}}\cdots\partial_{i_{s}}
\Big(\sqrt{\det(\mathbf{g})}\,g^{i_{1}j_{1}}\cdots g^{i_{s}j_{s}}\,H^{a}_{j_{1}\dots j_{s}}\Big)
=
\mathrm{(I)}+\mathrm{(II)}+\mathrm{(III)}.
\end{equation} where \[
\mathrm{(I)}
:=
\frac{1}{\sqrt{\det(\mathbf{g})}}\,
\partial_{i_{1}}\cdots\partial_{i_{s}}
\Big(g^{i_{1}j_{1}}\cdots g^{i_{s}j_{s}}\partial_{p}\big(\sqrt{\det(\mathbf{g})}\,g^{pq}G^{a}_{qj_{1}\dots j_{s}}\big)\Big),
\] \[
\mathrm{(II)}
:=
-\frac{1}{\sqrt{\det(\mathbf{g})}}\,
\partial_{i_{1}}\cdots\partial_{i_{s}}
\Big(\sqrt{\det(\mathbf{g})}\,g^{i_{1}j_{1}}\cdots g^{i_{s}j_{s}}
\displaystyle\sum_{\ell=1}^{s}\Gamma^{\alpha}_{p j_{\ell}}\,g^{pq}\,
G^{a}_{qj_{1}\dots j_{\ell-1}\alpha j_{\ell+1}\dots j_{s}}\Big),
\] \[
\mathrm{(III)}
:=
\frac{1}{\sqrt{\det(\mathbf{g})}}\,
\partial_{i_{1}}\cdots\partial_{i_{s}}
\Big(\sqrt{\det(\mathbf{g})}\,g^{i_{1}j_{1}}\cdots g^{i_{s}j_{s}}
\displaystyle\sum_{b=1}^{r}(A_{p})^{a}_{b}\,g^{pq}G^{b}_{qj_{1}\dots j_{s}}\Big).
\]

For the term $\mathrm{(I)}$, we write \[
\mathrm{(I)}
=
\frac{1}{\sqrt{\det(\mathbf{g})}}\,
\partial_{i_{1}}\cdots\partial_{i_{s}}
\Big(
\big(g^{i_{1}j_{1}}\cdots g^{i_{s}j_{s}}\big)\,
\partial_{p}\big(\sqrt{\det(\mathbf{g})}\,g^{pq}G^{a}_{qj_{1}\dots j_{s}}\big)
\Big).
\] We move the factor $g^{i_{1}j_{1}}\cdots g^{i_{s}j_{s}}$ inside $\partial_p$ using the identity \[
\big(g^{i_{1}j_{1}}\cdots g^{i_{s}j_{s}}\big)\,\partial_{p}(f)
=
\partial_{p}\Big(\big(g^{i_{1}j_{1}}\cdots g^{i_{s}j_{s}}\big)f\Big)
-
\partial_{p}\big(g^{i_{1}j_{1}}\cdots g^{i_{s}j_{s}}\big)\,f,
\qquad
\forall\,f\in C^{\infty}(U).
\] Applying it to $f=\sqrt{\det(\mathbf{g})}\,g^{pq}G^{a}_{qj_{1}\dots j_{s}}$ gives \begin{equation}\label{eq:I-principal-no-multi}
\begin{aligned}
\mathrm{(I)}
&=
\frac{1}{\sqrt{\det(\mathbf{g})}}\,
\partial_{i_{1}}\cdots\partial_{i_{s}}\partial_{p}
\Big(\sqrt{\det(\mathbf{g})}\,g^{pq}g^{i_{1}j_{1}}\cdots g^{i_{s}j_{s}}\,
G^{a}_{qj_{1}\dots j_{s}}\Big)\\
&\quad
-\frac{1}{\sqrt{\det(\mathbf{g})}}\,
\partial_{i_{1}}\cdots\partial_{i_{s}}
\Big(
\partial_{p}\big(g^{i_{1}j_{1}}\cdots g^{i_{s}j_{s}}\big)\,
\sqrt{\det(\mathbf{g})}\,g^{pq}G^{a}_{qj_{1}\dots j_{s}}
\Big).
\end{aligned}
\end{equation} The first term in \eqref{eq:I-principal-no-multi} is the leading term of order $s+1$. In the second term, the factor $\partial_p(g^{i_1j_1}\cdots g^{i_sj_s})$ already contains the derivative separated from the leading term. Only the $s$ outer derivatives $\partial_{i_1}\cdots\partial_{i_s}$ remain to be distributed between the coefficients and $\mathbf G$ by the Leibniz rule; even when they all act on $\mathbf G$, the order of differentiation of its components is at most $s$. Consequently, there are smooth coefficients $\bigl(B_{\beta}\bigr)^{a\,j_{1}\dots j_{s}q}_{c}$ with $|\beta|\le s$ such that \begin{equation}\label{eq:I-remainder}
\mathrm{(I)}
=
\frac{1}{\sqrt{\det(\mathbf{g})}}\,
\partial_{i_{1}}\cdots\partial_{i_{s}}\partial_{p}
\Big(\sqrt{\det(\mathbf{g})}\,g^{pq}g^{i_{1}j_{1}}\cdots g^{i_{s}j_{s}}\,
G^{a}_{qj_{1}\dots j_{s}}\Big)
+
\displaystyle\sum_{c=1}^{r}\displaystyle\sum_{|\beta|\le s}
\bigl(B_{\beta}\bigr)^{a\,j_{1}\dots j_{s}q}_{c}\,
\partial^{\beta}\!\bigl(G^{c}_{qj_{1}\dots j_{s}}\bigr).
\end{equation}

Now consider $\mathrm{(II)}$. Here no additional derivative $\partial_{p}$ acts on $\mathbf{G}$. Thus, expanding the operator $\partial_{i_{1}}\cdots\partial_{i_{s}}$ acting on the product \[
\sqrt{\det(\mathbf{g})}\,g^{i_{1}j_{1}}\cdots g^{i_{s}j_{s}}\,
\Gamma^{\alpha}_{p j_{\ell}}\,g^{pq}\,
G^{a}_{qj_{1}\dots j_{\ell-1}\alpha j_{\ell+1}\dots j_{s}},
\] by the Leibniz rule yields only terms in which at most $s$ derivatives act on the components of $\mathbf{G}$. Consequently, there are smooth coefficients $\bigl(\widetilde{B}_{\beta}\bigr)^{a\,j_{1}\dots j_{s}q}_{c}$ with $|\beta|\le s$ such that \begin{equation}\label{eq:II-bound}
\mathrm{(II)}
=
\displaystyle\sum_{c=1}^{r}\displaystyle\sum_{|\beta|\le s}
\bigl(\widetilde{B}_{\beta}\bigr)^{a\,j_{1}\dots j_{s}q}_{c}\,
\partial^{\beta}\!\bigl(G^{c}_{qj_{1}\dots j_{s}}\bigr).
\end{equation}

In $\mathrm{(III)}$ the dependence on $\mathbf{G}$ is linear and contains no additional derivatives of $\mathbf{G}$. By the Leibniz formula, each of the $s$ outer derivatives can act on that factor at most once more; thus only derivatives $\partial^\beta \mathbf{G}$ with $|\beta|\leq s$ arise. Therefore, there are smooth coefficients $\bigl(\widehat{B}_{\beta}\bigr)^{a\,j_{1}\dots j_{s}q}_{c}$ with $|\beta|\le s$ such that \begin{equation}\label{eq:III-bound}
\mathrm{(III)}
=
\displaystyle\sum_{c=1}^{r}\displaystyle\sum_{|\beta|\le s}
\bigl(\widehat{B}_{\beta}\bigr)^{a\,j_{1}\dots j_{s}q}_{c}\,
\partial^{\beta}\!\bigl(G^{c}_{qj_{1}\dots j_{s}}\bigr).
\end{equation}

Substituting \eqref{eq:I-remainder}, \eqref{eq:II-bound}, and \eqref{eq:III-bound} into \eqref{eq:split-I-II-III} gives \begin{equation}\label{eq:principal-plus-low-from-Hterm}
\begin{aligned}
&\frac{1}{\sqrt{\det(\mathbf{g})}}\,
\partial_{i_{1}}\cdots\partial_{i_{s}}
\Big(\sqrt{\det(\mathbf{g})}\,g^{i_{1}j_{1}}\cdots g^{i_{s}j_{s}}\,H^{a}_{j_{1}\dots j_{s}}\Big)\\
&\quad=\frac{1}{\sqrt{\det(\mathbf{g})}}\,
\partial_{p}\partial_{i_{1}}\cdots\partial_{i_{s}}
\Big(\sqrt{\det(\mathbf{g})}\,g^{pq}g^{i_{1}j_{1}}\cdots g^{i_{s}j_{s}}\,
G^{a}_{qj_{1}\dots j_{s}}\Big)\\
&\qquad+\displaystyle\sum_{c=1}^{r}\displaystyle\sum_{|\beta|\le s}
\bigl(D_{\beta}\bigr)^{a\,j_{1}\dots j_{s}q}_{c}\,
\partial^{\beta}\!\bigl(G^{c}_{qj_{1}\dots j_{s}}\bigr).
\end{aligned}
\end{equation} where \[
\bigl(D_{\beta}\bigr)^{a\,j_{1}\dots j_{s}q}_{c}
:=
\bigl(B_{\beta}\bigr)^{a\,j_{1}\dots j_{s}q}_{c}
+\bigl(\widetilde{B}_{\beta}\bigr)^{a\,j_{1}\dots j_{s}q}_{c}
+\bigl(\widehat{B}_{\beta}\bigr)^{a\,j_{1}\dots j_{s}q}_{c},
\qquad
|\beta|\le s.
\]

On the other hand, for the lower-order term in \eqref{eq:IH-on-H}, \[
\displaystyle\sum_{b=1}^{r}\displaystyle\sum_{|\beta|\le s-1}
\bigl(C_{\beta}\bigr)^{a\,j_{1}\dots j_{s}}_{b}\,
\partial^{\beta}\!\bigl(H^{b}_{j_{1}\dots j_{s}}\bigr),
\] substituting \eqref{eq:H-expanded} into $H^{b}_{j_{1}\dots j_{s}}$ shows that each $\partial^{\beta}(H^{b}_{j_{1}\dots j_{s}})$ is a linear combination (with smooth coefficients) of terms $\partial^{\gamma}(G^{c}_{qj_{1}\dots j_{s}})$ with $|\gamma|\le |\beta|+1\le s$. Consequently, there are smooth coefficients $\bigl(\overline{D}_{\beta}\bigr)^{a\,j_{1}\dots j_{s}q}_{c}$ with $|\beta|\le s$ such that \begin{equation}\label{eq:lowterm-from-IH}
\displaystyle\sum_{b=1}^{r}\displaystyle\sum_{|\beta|\le s-1}
\bigl(C_{\beta}\bigr)^{a\,j_{1}\dots j_{s}}_{b}\,
\partial^{\beta}\!\bigl(H^{b}_{j_{1}\dots j_{s}}\bigr)
=
\displaystyle\sum_{c=1}^{r}\displaystyle\sum_{|\beta|\le s}
\bigl(\overline{D}_{\beta}\bigr)^{a\,j_{1}\dots j_{s}q}_{c}\,
\partial^{\beta}\!\bigl(G^{c}_{qj_{1}\dots j_{s}}\bigr).
\end{equation}

Combining \eqref{eq:IH-on-H}, \eqref{eq:principal-plus-low-from-Hterm}, and \eqref{eq:lowterm-from-IH}, we conclude that there are smooth coefficients $\bigl(C_{\beta}\bigr)^{a\,j_{1}\dots j_{s}q}_{c}$ with $|\beta|\le s$ such that \[
\begin{aligned}
\bigl((\operatorname{tr}_{\mathbf{g}}\circ\nabla)^{s+1}\mathbf{G}\bigr)^{a}
&=\frac{1}{\sqrt{\det(\mathbf{g})}}\,
\partial_{i_{1}}\cdots\partial_{i_{s}}\partial_{p}
\Big(\sqrt{\det(\mathbf{g})}\,g^{pq}g^{i_{1}j_{1}}\cdots g^{i_{s}j_{s}}\,
G^{a}_{qj_{1}\dots j_{s}}\Big)\\
&\quad+\displaystyle\sum_{c=1}^{r}\displaystyle\sum_{|\beta|\le s}
\bigl(C_{\beta}\bigr)^{a\,j_{1}\dots j_{s}q}_{c}\,
\partial^{\beta}\!\bigl(G^{c}_{qj_{1}\dots j_{s}}\bigr).
\end{aligned}
\] Finally, we rename $p=i_1$ and $q=j_1$, simultaneously replacing the previous indices $(i_\ell,j_\ell)$ by $(i_{\ell+1},j_{\ell+1})$ for $\ell\in\{1,\dots,s\}$. The partial derivatives commute, so their outer order can be rearranged, while the indices of $\mathbf G$ remain in their original order. This gives exactly the form of \eqref{eq:claim} for $s+1$, completing the induction.

\end{proof} Let us consider some examples illustrating the preceding result. \begin{example}\label{ej:definicion-sobolev-haces-abierto-variedad-riemanniana-metrica-euclidiana-haz} Let $\Omega\subseteq\mathbb R^n$ be open, and regard $M:=\Omega$ as a smooth manifold without boundary, equipped with the Euclidean metric $\overline{\mathbf{g}}$. In the trivial line bundle \[
\mathbf{E}:=M\times\mathbb R\longrightarrow M
\], let $\mathbf{e}(p):=(p,1)$ be the constant unit section. Every section has a unique expression $\mathbf{u}=u\mathbf{e}$, with $u\in C^\infty(\Omega)$; the trivial metric and connection satisfy \[
\mathbf{h}_{\mathbf{E}}(u\mathbf{e},v\mathbf{e})=uv,
\qquad
\nabla^{\mathbf{E}}_{\boldsymbol{\partial}_i}(u\mathbf{e})
=(\partial_i u)\mathbf{e}.
\] For each $s\geq0$, the bundle isomorphism \[
I_s\colon T^{(0,s)}(TM)\otimes\mathbf{E}
\longrightarrow T^{(0,s)}(TM),
\qquad
I_s\bigl(\boldsymbol{\xi}\otimes\lambda\mathbf{e}(p)\bigr)
:=\lambda\boldsymbol{\xi},
\] identifies $\nabla^s\mathbf{u}$ with the Euclidean derivative $D^su$. In particular, \[
\|\mathbf{u}\|_{W^{m,p}(M,\mathbf{E})}
=
\left(\sum_{s=0}^m\|D^su\|_{L^p(\Omega)}^p\right)^{1/p}
\] if $1\leq p<\infty$, with the usual modification using maxima when $p=\infty$. The same identification in the weak derivative identity shows that \[
\mathbf{u}\in W^{m,p}(M,\mathbf{E})
\quad\Longleftrightarrow\quad
u\in W^{m,p}(\Omega).
\] Thus the bundle definition recovers precisely the scalar Sobolev space for the trivial rank-one bundle. \end{example} The formal adjoint allows us to define certain nonlinear operators important in the study of PDEs, such as the $p$--Laplacian, in the setting of vector bundles.

\begin{definition}[The scalar and vector bundle $p$-Laplacian]\label{def:definicion-sobolev-haces-laplaciano-escalar-y-en-haces-vectoriales}\index{p Laplacian scalar and vector bundle@p-Laplacian, scalar and vector bundle} Let $(M,\mathbf{g})$ be a Riemannian manifold without boundary, and let $p>1$.

\begin{enumerate}[label=(\alph*)]

\item If $u\in C^{2}(M)$, the \emph{scalar $p$-Laplacian} of $u$ is defined by \[
\Delta_{p}u
:=
\operatorname{div}\bigl(|\nabla u|_{\mathbf{g}}^{p-2}\nabla u\bigr).
\]

Equivalently, using the formal adjoint of the covariant derivative, \[
\Delta_{p}u
=
-\nabla^{*}\bigl(|\nabla u|_{\mathbf{g}}^{p-2}\nabla u\bigr).
\]

When $p=2$, this recovers the Laplace--Beltrami Laplacian $\Delta_{2}u=\Delta_{\mathbf{g}}u$.

\item Now let $\mathbf{E}\to M$ be a smooth vector bundle equipped with a bundle metric $\mathbf{h}_{\mathbf{E}}$ and a compatible connection $\nabla^{\mathbf{E}}$. If $\mathbf{u}\in \Gamma^{2}(\mathbf{E})$, the \emph{$p$-Laplacian associated with $\nabla^{\mathbf{E}}$} is defined by \[
\Delta_{p}^{\mathbf{E}}\mathbf{u}
:=
-\left(\nabla^{\mathbf{E}}\right)^{*}
\Bigl(
|\nabla^{\mathbf{E}}\mathbf{u}|_{\mathbf{g},\mathbf{h}_{\mathbf{E}}}^{p-2}\nabla^{\mathbf{E}}\mathbf{u}
\Bigr).
\]

\end{enumerate} In both definitions we adopt $|\xi|^{p-2}\xi=0$ when $\xi=0$. If $1<p<2$, even for data of class $C^2$, this flux may fail to be of class $C^1$ at its zeros. In that case, the divergence and adjoint in the preceding formulas are interpreted weakly: for a test section $\boldsymbol\varphi$ supported in the interior, the value of $\Delta_p^{\mathbf E}\mathbf u$ is determined by \[
\big\langle\Delta_p^{\mathbf E}\mathbf u,\boldsymbol\varphi\big\rangle
=-\int_M\left\langle|\nabla^{\mathbf E}\mathbf u|^{p-2}\nabla^{\mathbf E}\mathbf u,\nabla^{\mathbf E}\boldsymbol\varphi\right\rangle_{\mathbf g,\mathbf h_{\mathbf E}}\,d\lambda_{\mathbf g}.
\] The scalar case is obtained using the trivial bundle. The flux is continuous and hence locally integrable, so this identity is always defined. When the flux is of class $C^1$, it recovers the classical expression, and the following boundary formula applies. For $p\geq2$, data of class $C^2$ ensure this regularity of the flux. \end{definition} Our integration-by-parts formula also yields further integration-by-parts formulas, such as the following: \begin{proposition}[Integration-by-parts formula for the \(p\)-Laplacian]\label{prop:definicion-sobolev-haces-formula-de-integracion-por-partes-para-el-p-laplaciano}\index{integration by parts formula for the p Laplacian@integration-by-parts formula for the \(p\)-Laplacian} Let $(M,\mathbf{g})$ be a Riemannian manifold with or without boundary, and let $p>1$.

In the scalar case, let $u\in C^{2}(M)$ be such that $|\nabla u|_{\mathbf{g}}^{p-2}\nabla u$ is of class $C^{1}$. Then, for every $\varphi\in C^{\infty}_{c}(M)$, \[
\int_{M}|\nabla u|_{\mathbf{g}}^{p-2}\langle \nabla u,\nabla\varphi\rangle_{\mathbf{g}}\,d\lambda_{\mathbf{g}}
=
-\int_{M}\varphi\,\Delta_{p}u\,d\lambda_{\mathbf{g}}
+
\int_{\partial M}\varphi\,|\nabla u|_{\mathbf{g}}^{p-2}\partial_{\boldsymbol{\nu}}u\,d\lambda_{\widetilde{\mathbf{g}}}.
\]

More generally, let $\pi_{\mathbf{E}}\colon \mathbf{E}\longrightarrow M$ be a smooth vector bundle equipped with a bundle metric $\mathbf{h}_{\mathbf{E}}$ and a compatible connection $\nabla^{\mathbf{E}}$. If $\mathbf{u}\in\Gamma^{2}(\mathbf{E})$ is such that $|\nabla^{\mathbf{E}}\mathbf{u}|_{\mathbf{g},\mathbf{h}_{\mathbf{E}}}^{p-2}\nabla^{\mathbf{E}}\mathbf{u}$ is of class $C^{1}$, then, for every $\boldsymbol{\phi}\in\Gamma_{c}(\mathbf{E})$, \[
\int_{M}
|\nabla^{\mathbf{E}}\mathbf{u}|_{\mathbf{g},\mathbf{h}_{\mathbf{E}}}^{p-2}
\langle \nabla^{\mathbf{E}}\mathbf{u},\nabla^{\mathbf{E}}\boldsymbol{\phi}\rangle_{\mathbf{g},\mathbf{h}_{\mathbf{E}}}\,d\lambda_{\mathbf{g}}
=
-\int_{M}
\langle \Delta_{p}^{\mathbf{E}}\mathbf{u},\boldsymbol{\phi}\rangle_{\mathbf{h}_{\mathbf{E}}}\,d\lambda_{\mathbf{g}}
+
\int_{\partial M}
\left\langle
\iota_{\boldsymbol{\nu}}\bigl(|\nabla^{\mathbf{E}}\mathbf{u}|_{\mathbf{g},\mathbf{h}_{\mathbf{E}}}^{p-2}\nabla^{\mathbf{E}}\mathbf{u}\bigr),\boldsymbol{\phi}
\right\rangle_{\mathbf{h}_{\mathbf{E}}}\,d\lambda_{\widetilde{\mathbf{g}}}.
\] In particular, if $\boldsymbol{\phi}$ has compact support contained in $\operatorname{Int}(M)$, the boundary term vanishes. \end{proposition}

\begin{proof} We first prove the bundle identity. Define $\mathbf{A}_{\mathbf{u}}:=|\nabla^{\mathbf{E}}\mathbf{u}|_{\mathbf{g},\mathbf{h}_{\mathbf{E}}}^{p-2}\nabla^{\mathbf{E}}\mathbf{u}$, so that $\mathbf{A}_{\mathbf{u}}\in\Gamma^{1}(T^{*}M\otimes \mathbf{E})$ by the regularity hypothesis. The calculation in the case $s=1$ of Green's formula uses only continuous first derivatives, so it remains valid with this regularity. We also multiply $\mathbf A_{\mathbf u}$ by a compactly supported smooth function equal to one on a neighborhood of $\operatorname{supp}\boldsymbol\phi$; this leaves every integral unchanged and satisfies the support condition. We apply Theorem~\ref{teo:green-m-esima-covariante-E} with $s=1$, $k=l=0$, $\mathbf{F}=\boldsymbol{\phi}$, and $\mathbf{G}=\mathbf{A}_{\mathbf{u}}$. Then \[
\int_{M}\langle \nabla^{\mathbf{E}}\boldsymbol{\phi},\mathbf{A}_{\mathbf{u}}\rangle_{\mathbf{g},\mathbf{h}_{\mathbf{E}}}\,d\lambda_{\mathbf{g}}
=
-\int_{M}
\big\langle \boldsymbol{\phi},(\operatorname{tr}_{\mathbf{g}}\circ\nabla)\mathbf{A}_{\mathbf{u}}\big\rangle_{\mathbf{h}_{\mathbf{E}}}\,d\lambda_{\mathbf{g}}
+
\int_{\partial M}
\langle \boldsymbol{\phi},\iota_{\boldsymbol{\nu}}\mathbf{A}_{\mathbf{u}}\rangle_{\mathbf{h}_{\mathbf{E}}}\,d\lambda_{\widetilde{\mathbf{g}}}.
\] Since $(\nabla^{\mathbf{E}})^{*}\mathbf{A}_{\mathbf{u}}=-(\operatorname{tr}_{\mathbf{g}}\circ\nabla)\mathbf{A}_{\mathbf{u}}$ and, by definition, $\Delta_{p}^{\mathbf{E}}\mathbf{u}=-(\nabla^{\mathbf{E}})^{*}\mathbf{A}_{\mathbf{u}}$, we conjugate the preceding identity in the complex case, or use symmetry in the real case. The interior term then becomes $-\langle \Delta_{p}^{\mathbf{E}}\mathbf{u},\boldsymbol{\phi}\rangle_{\mathbf{h}_{\mathbf{E}}}$. Moreover, with the bar acting trivially in the real case, $\overline{\langle \nabla^{\mathbf{E}}\boldsymbol{\phi},\mathbf{A}_{\mathbf{u}}\rangle_{\mathbf{g},\mathbf{h}_{\mathbf{E}}}}
=
|\nabla^{\mathbf{E}}\mathbf{u}|_{\mathbf{g},\mathbf{h}_{\mathbf{E}}}^{p-2}
\langle \nabla^{\mathbf{E}}\mathbf{u},\nabla^{\mathbf{E}}\boldsymbol{\phi}\rangle_{\mathbf{g},\mathbf{h}_{\mathbf{E}}}$. Substitution of these two identities gives the bundle formula.

The scalar case follows by taking $\mathbf{E}=M\times\mathbb R$ with its trivial metric and connection and identifying a section with its component function $u$. In this case $\mathbf{A}_{u}=|\nabla u|_{\mathbf{g}}^{p-2}\nabla u$ and \[
\iota_{\boldsymbol{\nu}}\mathbf{A}_{u}
=
\langle \mathbf{A}_{u},\boldsymbol{\nu}\rangle_{\mathbf{g}}
=
|\nabla u|_{\mathbf{g}}^{p-2}\langle\nabla u,\boldsymbol{\nu}\rangle_{\mathbf{g}}
=
|\nabla u|_{\mathbf{g}}^{p-2}\partial_{\boldsymbol{\nu}}u.
\] Thus the preceding identity becomes exactly the scalar formula. \end{proof}

The integration-by-parts formula also allows us to compare the usual covariant norm on $H^m(M,\mathbf E)$ with a norm defined by powers of the Laplacian. This equivalence is a recurring tool in geometric analysis, since it translates control of derivatives into spectral control. The next section establishes the result with the hypotheses and constants needed in the vector bundle setting.

\subsection{Equivalence of norms: an application of the integration-by-parts formula in vector bundles} In the Euclidean case, the integration-by-parts formula is fundamental to Sobolev space theory and nonlinear analysis. In this section we use it to prove, on a closed manifold $M$, the equivalence of the usual norm on $H^{m}(M):=W^{m,2}(M)$ and a norm defined by powers of the Laplacian. We begin with the corresponding identity on Euclidean domains: \begin{theorem}[An equivalent norm on $H_0^m(\Omega)$ in terms of the Laplacian] \label{teo:sobolev_abierto_euclidiano}\index{Sobolev space@Sobolev space!equivalent norm using the Laplacian} Let $\Omega \subset \mathbb{R}^n$ be an open set, and let $m \in \mathbb{N}$. For every function $u \in H_0^m(\Omega)$, the standard Sobolev norm of order $m$, defined using the partial derivative operator, \[
||u||_{H^m(\Omega)}^2 := \displaystyle\sum_{j=0}^m \int_{\Omega} |\nabla^j u|^2 \, d\lambda_n,
\] is exactly equal to the norm structured in powers of the Laplacian operator $\Delta = \displaystyle\sum_{i=1}^n \partial_i^2$, given by \[
||u||_{\Delta}^2 := \displaystyle\sum_{\substack{j=0 \\ j \text{ even}}}^m \int_{\Omega} |\Delta^{\frac{j}{2}} u|^2 \, d\lambda_n + \displaystyle\sum_{\substack{j=1 \\ j \text{ odd}}}^m \int_{\Omega} \big|\nabla \Delta^{\frac{j-1}{2}} u\big|^2 \, d\lambda_n.
\] Moreover, equality holds term by term for each order of differentiation $j \in \{0, \dots, m\}$. \end{theorem}

\begin{proof} First let $u\in H_0^m(\Omega)$. By the definition of this space, there is a sequence $(u_\nu)_{\nu\in\mathbb N}\subseteq C_c^\infty(\Omega)$ such that $u_\nu\to u$ in $H^m(\Omega)$. In particular, for every multi-index $\alpha$ with $|\alpha|\leq m$, \[
D^\alpha u_\nu\longrightarrow D^\alpha u
\qquad\text{in }L^2(\Omega).
\] If $j=2k\leq m$, the operator $\Delta^k$ is a finite linear combination, with constant coefficients, of partial derivatives of order $j$. If $j=2k+1\leq m$, each component of $\nabla\Delta^k$ is likewise a finite linear combination of partial derivatives of order $j$. Therefore, \[
\Delta^k u_\nu\longrightarrow\Delta^k u
\quad\text{and}\quad
\nabla\Delta^k u_\nu\longrightarrow\nabla\Delta^k u
\] in $L^2(\Omega)$ whenever the corresponding order does not exceed $m$. Thus, once the term-by-term equalities have been proved for $u_\nu$, we can pass to the limit in each of them. It therefore suffices to establish them for a function $u\in C_c^\infty(\Omega)$. For such a function, all integrations by parts have no boundary terms. Moreover, partial derivatives commute, and all norms are taken with respect to the Euclidean metric, by contraction of indices.

We proceed by induction on $j$ and prove, for every $u\in C_c^\infty(\Omega)$, the equality corresponding to order $j$.

For $j=0$ and $j=1$ the assertion follows immediately from the definition.

Now consider the fundamental case $j=2$. Observe that \[
\int_\Omega |\nabla^2 u|^2\,d\lambda_n
=
\int_\Omega \displaystyle\sum_{p,q=1}^n (\partial_p\partial_q u)^2\,d\lambda_n.
\] We first integrate by parts with respect to the variable $x^q$ and then with respect to $x^p$. Since $u$ has compact support in $\Omega$ and partial derivatives commute, we obtain \[
\int_\Omega (\partial_p\partial_q u)^2\,d\lambda_n
=
-\int_\Omega (\partial_p u)\,\partial_p\partial_q^2 u\,d\lambda_n
=
\int_\Omega (\partial_p^2 u)\,(\partial_q^2 u)\,d\lambda_n.
\] Summing over $p,q$ gives \[
\int_\Omega |\nabla^2 u|^2\,d\lambda_n
=
\int_\Omega \left(\displaystyle\sum_{p=1}^n \partial_p^2 u\right)
\left(\displaystyle\sum_{q=1}^n \partial_q^2 u\right)\,d\lambda_n
=
\int_\Omega |\Delta u|^2\,d\lambda_n.
\]

Now suppose that the assertion holds for all orders less than $j\geq 3$. For each list $(i_1,\dots,i_{j-2})$, the function $\partial_{i_1}\cdots\partial_{i_{j-2}}u$ belongs to $C_c^\infty(\Omega)$. Applying the identity already proved for two derivatives to each of these functions, and summing over all indices, gives, in tensor notation, \[
\int_\Omega |\nabla^j u|^2\,d\lambda_n
=
\int_\Omega |\Delta(\nabla^{j-2}u)|^2\,d\lambda_n.
\] Since the Laplacian commutes with partial derivatives in $\mathbb{R}^n$, we have \[
\Delta(\nabla^{j-2}u)=\nabla^{j-2}(\Delta u),
\] and therefore \[
\int_\Omega |\nabla^j u|^2\,d\lambda_n
=
\int_\Omega |\nabla^{j-2}(\Delta u)|^2\,d\lambda_n.
\]

The function $v:=\Delta u$ also belongs to $C_c^\infty(\Omega)$. We can therefore apply the induction hypothesis to it at order $j-2$.

If $j$ is even, then $j-2$ is even, and we obtain \[
\int_\Omega |\nabla^{j-2}(\Delta u)|^2
=
\int_\Omega |\Delta^{\frac{j-2}{2}}(\Delta u)|^2
=
\int_\Omega |\Delta^{\frac{j}{2}}u|^2.
\]

If $j$ is odd, then $j-2$ is odd, and we obtain \[
\int_\Omega |\nabla^{j-2}(\Delta u)|^2
=
\int_\Omega \big|\nabla \Delta^{\frac{j-3}{2}}(\Delta u)\big|^2
=
\int_\Omega \big|\nabla \Delta^{\frac{j-1}{2}}u\big|^2.
\]

In both cases this recovers exactly the corresponding term in the structured norm. This proves the term-by-term identities for functions in $C_c^\infty(\Omega)$. Applying them to the sequence $(u_\nu)$ chosen at the beginning and using the convergence in $L^2(\Omega)$ established there, we pass to the limit on both sides of each identity. This gives the same equalities for every $u\in H_0^m(\Omega)$; summing them for $0\leq j\leq m$ yields equality of the two norms. \end{proof}

To extend this result to Riemannian manifolds, we first introduce notation that simplifies the manipulation of tensor contractions, together with an auxiliary analytic result. \begin{notation}\label{notacion estrella} Let $(M,\mathbf{g})$ be a Riemannian manifold with or without boundary, and let \[
\pi_{\mathbf{E}}\colon \mathbf{E}\longrightarrow M,
\qquad
\pi_{\mathbf{F}}\colon \mathbf{F}\longrightarrow M
\] be smooth vector bundles of finite rank, equipped with smooth bundle metrics $\mathbf{h}_{\mathbf{E}}$ and $\mathbf{h}_{\mathbf{F}}$.

Let $\mathbf{G}_1,\dots,\mathbf{G}_r,\mathbf{H}$ be vector bundles obtained from $\mathbf{E}$, $\mathbf{F}$, their duals, and the tangent and cotangent bundles of $M$ by finitely many tensor products. Let also \[
\mathbf{u}_1\in \Gamma(\mathbf{G}_1),\dots,\mathbf{u}_r\in \Gamma(\mathbf{G}_r).
\]

We say that a section $\mathbf{T}\in \Gamma(\mathbf{H})$ has the form \[
\mathbf{u}_1 * \cdots * \mathbf{u}_r
\] if there is a smooth bundle homomorphism \[
\Phi\colon \mathbf{G}_1 \otimes \cdots \otimes \mathbf{G}_r \longrightarrow \mathbf{H}
\] obtained by composition and finite linear combination of the following natural operations:

\begin{enumerate}  \item tensor products;  \item permutations of tensor factors;  \item contractions between a bundle and its dual;  \item the bundle isomorphisms induced by the bundle metrics $\mathbf{h}_{\mathbf{E}}$, $\mathbf{h}_{\mathbf{F}}$ and the Riemannian metric $\mathbf{g}$, together with their inverses. \end{enumerate}

and such that \[
\mathbf{T} = \Phi(\mathbf{u}_1 \otimes \cdots \otimes \mathbf{u}_r).
\]

In particular, if $\mathbf{u}\in \Gamma(\mathbf{G})$ and $\mathbf{v}\in \Gamma(\mathbf{K})$, the expression \[
\mathbf{u} * \mathbf{v}
\] denotes any section obtained from $\mathbf{u}\otimes \mathbf{v}$ by a finite combination of the preceding operations. Different occurrences of the symbol $*$ may correspond to different bundle homomorphisms. \end{notation}

Let us consider an example illustrating how this notation works.

\begin{example}\label{ej:definicion-sobolev-haces-variedad-riemanniana-isomorfismo-canonico-haces-homomorfismo} Let $(M,\mathbf{g})$ be a Riemannian manifold with or without boundary. Consider \[
\mathbf{R} \in \Gamma(T^{(0,4)}(TM)),
\qquad
\boldsymbol{\omega} \in \Gamma(T^{(0,1)}(TM)).
\]

There is a canonical bundle isomorphism \[
T^{(0,4)}(TM)\otimes T^{(0,1)}(TM)
\cong
T^{(0,5)}(TM).
\]

Define a bundle homomorphism \[
\Phi\colon T^{(0,5)}(TM)\longrightarrow T^{(0,1)}(TM)
\] by metric contractions with $\mathbf{g}^{-1}\in \Gamma(T^{(2,0)}(TM))$. More explicitly, in local coordinates, if \[
\mathbf{A}\in \Gamma(T^{(0,5)}(TM)),
\] and $A_{abcde}$ are its components, then the components of $\Phi(\mathbf{A})$ are \[
\Phi(\mathbf{A})_b = g^{ae} g^{cd} A_{abcde}.
\]

This homomorphism is obtained from tensor products, permutations, and contractions, so it is among those allowed in Notation \ref{notacion estrella}.

Evaluating at $\mathbf{A} = \mathbf{R}\otimes \boldsymbol{\omega}$ gives \[
\Phi(\mathbf{R}\otimes \boldsymbol{\omega})_i
=
g^{jk} g^{lm} R_{jilm}\,\omega_k.
\]

Consequently, \[
\mathbf{T} = \Phi(\mathbf{R}\otimes \boldsymbol{\omega})
\] is a term of the type $\mathbf{R} * \boldsymbol{\omega}$. \end{example} \begin{lemma}\label{lema estimacion estrella} Let $(M,\mathbf{g})$ be a Riemannian manifold with or without boundary, and let \[
\mathbf{E}_1,\dots,\mathbf{E}_r,\mathbf{H} \longrightarrow M
\] be smooth vector bundles of finite rank, equipped with smooth bundle metrics.

Let \[
\Phi\colon \mathbf{E}_1\otimes\cdots\otimes \mathbf{E}_r\longrightarrow \mathbf{H}
\] be a smooth bundle homomorphism obtained by composition and finite linear combination of tensor products, permutations of factors, contractions, and the isomorphisms induced by the bundle metrics and the Riemannian metric $\mathbf{g}$, together with their inverses.

Then, for every compact set $K\subseteq M$, there is a constant $C_{K}>0$ such that, for every $x\in K$ and all \[
v_i\in (\mathbf{E}_i)_x,\qquad i=1,\dots,r,
\], we have \[
\bigl|\Phi_x(v_1\otimes\cdots\otimes v_r)\bigr|_{\mathbf{h}_{\mathbf{H}}}
\leq
C_{K}\,|v_1|_{\mathbf{h}_{\mathbf{E}_1}}\cdots |v_r|_{\mathbf{h}_{\mathbf{E}_r}}.
\]

In particular, if $\mathbf{u}_i\in \Gamma(\mathbf{E}_i)$ for $i=1,\dots,r$, then \[
|\Phi(\mathbf{u}_1\otimes\cdots\otimes \mathbf{u}_r)(x)|_{\mathbf{h}_{\mathbf{H}}}
\leq
C_{K}\,|\mathbf{u}_1(x)|_{\mathbf{h}_{\mathbf{E}_1}}\cdots |\mathbf{u}_r(x)|_{\mathbf{h}_{\mathbf{E}_r}},
\qquad \forall x\in K.
\] Consequently, if $T=\mathbf{u}_1*\cdots *\mathbf{u}_r$, then \[
|T(x)|_{\mathbf{h}_{\mathbf{H}}}\leq C_{K}\,|\mathbf{u}_1(x)|_{\mathbf{h}_{\mathbf{E}_1}}\cdots |\mathbf{u}_r(x)|_{\mathbf{h}_{\mathbf{E}_r}},
\qquad \forall x\in K.
\] \end{lemma}

\begin{proof} Fix \(x\in M\). Consider the map \[
\widetilde{\Phi}_x\colon (\mathbf{E}_1)_x\times\cdots\times (\mathbf{E}_r)_x\longrightarrow \mathbf{H}_x,
\qquad
\widetilde{\Phi}_x(v_1,\dots,v_r):=\Phi_x(v_1\otimes\cdots\otimes v_r).
\]

Let \[
\theta_x:(\mathbf{E}_1)_x\times\cdots\times (\mathbf{E}_r)_x
\longrightarrow
(\mathbf{E}_1)_x\otimes\cdots\otimes (\mathbf{E}_r)_x,
\qquad
\theta_x(v_1,\dots,v_r)=v_1\otimes\cdots\otimes v_r.
\] The map \(\theta_x\) is multilinear. On the other hand, since \(\Phi\) is a bundle homomorphism, \[
\Phi_x\colon (\mathbf{E}_1)_x\otimes\cdots\otimes (\mathbf{E}_r)_x\longrightarrow \mathbf{H}_x
\] is linear. Consequently, \[
\widetilde{\Phi}_x=\Phi_x\circ \theta_x
\] is multilinear.

Denote by \(|\cdot|_{\mathbf{h}_{\mathbf{E}_i,x}}\) the norm on \((\mathbf{E}_i)_x\) induced by the bundle metric \(\mathbf{h}_{\mathbf{E}_i}\), and by \(|\cdot|_{\mathbf{h}_{\mathbf{H},x}}\) the norm on \(\mathbf{H}_x\) induced by the bundle metric of \(\mathbf{H}\). Since the spaces involved are finite-dimensional, every multilinear map between them is continuous. Thus there is a constant \(C(x)>0\) such that \[
|\Phi_x(v_1\otimes\cdots\otimes v_r)|_{\mathbf{h}_{\mathbf{H},x}}
\le
C(x)\,
|v_1|_{\mathbf{h}_{\mathbf{E}_1,x}}\cdots |v_r|_{\mathbf{h}_{\mathbf{E}_r,x}}
\] for any \(v_i\in (\mathbf{E}_i)_x\), with \(i\in\{1,\dots,r\}\).

It remains to show that this constant can be chosen uniformly on compact sets. Let \(p\in K\). Choose open sets \(U_p,V_p\subseteq M\) such that \[
p\in U_p,
\qquad
\overline{U_p}\subseteq V_p,
\qquad
\overline{U_p}\text{ is compact},
\] and such that all the bundles \[
\mathbf{E}_1,\dots,\mathbf{E}_r,\mathbf{H}
\] admit smooth local frames over \(V_p\). Fix local frames on \(V_p\), \[
e^{(i)}_1,\dots,e^{(i)}_{m_i}
\quad \text{for } \mathbf{E}_i,
\qquad
\boldsymbol{\varepsilon}_1,\dots,\boldsymbol{\varepsilon}_{\ell}
\quad \text{for } \mathbf{H}.
\]

In these frames there are smooth functions \[
\Phi^b_{a_1\cdots a_r}\colon V_p\longrightarrow \mathbb{R}
\] such that \[
\Phi_x\bigl(
e^{(1)}_{a_1}(x)\otimes\cdots\otimes e^{(r)}_{a_r}(x)
\bigr)
=
\displaystyle\sum_{b=1}^{\ell}
\Phi^b_{a_1\cdots a_r}(x)\,\boldsymbol{\varepsilon}_b(x).
\]

Since the frames and these functions are defined on a neighborhood of \(\overline{U_p}\), their continuity and the compactness of \(\overline{U_p}\) give \(A_p>0\) such that \[
|\Phi^b_{a_1\cdots a_r}(x)|\le A_p
\qquad \forall x\in \overline{U_p}.
\]

Now let \(x\in \overline{U_p}\) and, for each \(i\in\{1,\dots,r\}\), write \[
v_i=\displaystyle\sum_{a_i=1}^{m_i} v_i^{a_i}\,e^{(i)}_{a_i}(x).
\] Then \[
\Phi_x(v_1\otimes\cdots\otimes v_r)
=
\displaystyle\sum_{\substack{
1\le b\le \ell\\
1\le a_i\le m_i
}}
\Phi^b_{a_1\cdots a_r}(x)\,
v_1^{a_1}\cdots v_r^{a_r}\,
\boldsymbol{\varepsilon}_b(x).
\]

By the triangle inequality, \[
|\Phi_x(v_1\otimes\cdots\otimes v_r)|_{\mathbf{h}_{\mathbf{H},x}}
\le
\displaystyle\sum_{\substack{
1\le b\le \ell\\
1\le a_i\le m_i
}}
|\Phi^b_{a_1\cdots a_r}(x)|\,
|v_1^{a_1}|\cdots |v_r^{a_r}|\,
|\boldsymbol{\varepsilon}_b(x)|_{\mathbf{h}_{\mathbf{H},x}}.
\] Using the bound \(|\Phi^b_{a_1\cdots a_r}(x)|\le A_p\) gives \[
|\Phi_x(v_1\otimes\cdots\otimes v_r)|_{\mathbf{h}_{\mathbf{H},x}}
\le
A_p
\displaystyle\sum_{\substack{
1\le b\le \ell\\
1\le a_i\le m_i
}}
|v_1^{a_1}|\cdots |v_r^{a_r}|\,
|\boldsymbol{\varepsilon}_b(x)|_{\mathbf{h}_{\mathbf{H},x}}.
\]

Since \(\overline{U_p}\) is compact and the sections \(\boldsymbol{\varepsilon}_b\) are smooth, there is \(B_p>0\) such that \[
|\boldsymbol{\varepsilon}_b(x)|_{\mathbf{h}_{\mathbf{H},x}}\le B_p
\qquad \forall x\in \overline{U_p},\ \forall b.
\] Therefore, \[
|\Phi_x(v_1\otimes\cdots\otimes v_r)|_{\mathbf{h}_{\mathbf{H},x}}
\le
A_pB_p
\displaystyle\sum_{\substack{
1\le b\le \ell\\
1\le a_i\le m_i
}}
|v_1^{a_1}|\cdots |v_r^{a_r}|.
\]

Now, \[
\displaystyle\sum_{\substack{
1\le b\le \ell\\
1\le a_i\le m_i
}}
|v_1^{a_1}|\cdots |v_r^{a_r}|
=
\displaystyle\sum_{b=1}^{\ell}
\displaystyle\sum_{a_1=1}^{m_1}\cdots \displaystyle\sum_{a_r=1}^{m_r}
|v_1^{a_1}|\cdots |v_r^{a_r}|.
\] Since the inner expression is independent of \(b\), we obtain \[
\displaystyle\sum_{\substack{
1\le b\le \ell\\
1\le a_i\le m_i
}}
|v_1^{a_1}|\cdots |v_r^{a_r}|
=
\ell
\displaystyle\sum_{a_1=1}^{m_1}\cdots \displaystyle\sum_{a_r=1}^{m_r}
|v_1^{a_1}|\cdots |v_r^{a_r}|.
\] Moreover, by distributivity, \[
\displaystyle\sum_{a_1=1}^{m_1}\cdots \displaystyle\sum_{a_r=1}^{m_r}
|v_1^{a_1}|\cdots |v_r^{a_r}|
=
\left(\displaystyle\sum_{a_1=1}^{m_1}|v_1^{a_1}|\right)\cdots
\left(\displaystyle\sum_{a_r=1}^{m_r}|v_r^{a_r}|\right).
\] Consequently, \[
|\Phi_x(v_1\otimes\cdots\otimes v_r)|_{\mathbf{h}_{\mathbf{H},x}}
\le
A_pB_p\,\ell
\prod_{i=1}^{r}
\left(\displaystyle\sum_{a_i=1}^{m_i}|v_i^{a_i}|\right).
\]

For each \(i\in\{1,\dots,r\}
 \), the local frame \(\mathbf{e}^{(i)}_1,\dots,\mathbf{e}^{(i)}_{m_i}\) induces on each fiber \((\mathbf{E}_i)_x\), with \(x\in \overline{U_p}\), the norm \[
\|v_i\|_{(i),x}:=\displaystyle\sum_{a_i=1}^{m_i}|v_i^{a_i}|.
\] The matrix of $\mathbf{h}_{\mathbf{E}_i}$ in this frame and its inverse depend smoothly on $x$. Their coefficients are bounded on $\overline{U_p}$, and their positive minors are bounded away from zero. The compactness of $\overline{U_p}$ and of the unit sphere in each component space gives constants \(0<d_{p,i}\leq D_{p,i}<+\infty\), independent of $x$, such that \[
 d_{p,i}|v_i|_{\mathbf{h}_{\mathbf{E}_i,x}}
 \leq\|v_i\|_{(i),x}
 \leq D_{p,i}\,|v_i|_{\mathbf{h}_{\mathbf{E}_i,x}}
\qquad \forall x\in \overline{U_p},\ \forall v_i\in (\mathbf{E}_i)_x.
\] Substituting these estimates into the preceding inequality gives \[
|\Phi_x(v_1\otimes\cdots\otimes v_r)|_{\mathbf{h}_{\mathbf{H},x}}
\le
A_pB_p\,\ell\,D_{p,1}\cdots D_{p,r}\,
|v_1|_{\mathbf{h}_{\mathbf{E}_1,x}}\cdots |v_r|_{\mathbf{h}_{\mathbf{E}_r,x}}
\] for every \(x\in \overline{U_p}\). Defining \[
C_p:=A_pB_p\,\ell\,D_{p,1}\cdots D_{p,r},
\] yields \[
|\Phi_x(v_1\otimes\cdots\otimes v_r)|_{\mathbf{h}_{\mathbf{H},x}}
\le
C_p\,
|v_1|_{\mathbf{h}_{\mathbf{E}_1,x}}\cdots |v_r|_{\mathbf{h}_{\mathbf{E}_r,x}}
\qquad \forall x\in \overline{U_p}.
\]

The open sets \(U_p\), with \(p\in K\), cover \(K\). By compactness there is a finite subcover \[
U_{p_1},\dots,U_{p_N}.
\] Defining \[
C_K:=\max\{C_{p_1},\dots,C_{p_N}\},
\] gives \[
|\Phi_x(v_1\otimes\cdots\otimes v_r)|_{\mathbf{h}_{\mathbf{H},x}}
\le
C_K\,
|v_1|_{\mathbf{h}_{\mathbf{E}_1,x}}\cdots |v_r|_{\mathbf{h}_{\mathbf{E}_r,x}}
\qquad \forall x\in K.
\]

The assertion for sections follows by taking \(v_i=\mathbf{u}_i(x)\). \end{proof} \begin{proposition}[Leibniz rule for the operator $*$]\label{prop: leibniz operador *}\index{Leibniz rule for the operator@Leibniz rule for the operator} Let $(M,\mathbf{g})$ be a Riemannian manifold with or without boundary. Let $\mathbf{E}_1,\dots,\mathbf{E}_r,\mathbf{F}\longrightarrow M$ be smooth vector bundles with bundle metrics (Hermitian in the complex case) and connections compatible with them. We use the Levi--Civita connection on the tangent and cotangent factors, and the induced connections on all products, duals, and homomorphism bundles. Let \[
\boldsymbol{\Phi}\colon \mathbf{E}_1\otimes\cdots\otimes \mathbf{E}_r\longrightarrow \mathbf{F}
\] be a smooth bundle homomorphism obtained using tensor products, permutations, contractions, and the isomorphisms induced by the bundle metrics and the Riemannian metric.

If $\mathbf{u}_i\in \Gamma(\mathbf{E}_i)$ for $i=1,\dots,r$, define \[
\mathbf{T}=\boldsymbol{\Phi}(\mathbf{u}_1\otimes\cdots\otimes \mathbf{u}_r)=\mathbf{u}_1*\cdots *\mathbf{u}_r.
\]

Then, for every $\mathbf{X}\in \mathfrak{X}(M)$, \[
\nabla_{\mathbf{X}} \mathbf{T}
=
\displaystyle\sum_{i=1}^{r}
\mathbf{u}_1*\cdots *(\nabla_{\mathbf{X}} \mathbf{u}_i)*\cdots *\mathbf{u}_r.
\] \end{proposition}

\begin{proof} By definition, \(\mathbf{T}=\boldsymbol{\Phi}(\mathbf{u}_1\otimes\cdots\otimes \mathbf{u}_r)\). We regard \(\boldsymbol{\Phi}\) as a section of \(\operatorname{Hom}(\mathbf{E}_1\otimes\cdots\otimes \mathbf{E}_r,\mathbf{F})\), with the connection induced by those on \(\mathbf{E}_1,\dots,\mathbf{E}_r\) and \(\mathbf{F}\).

By the definition of this connection, for every \(\mathbf{S}\in \Gamma(\mathbf{E}_1\otimes\cdots\otimes \mathbf{E}_r)\) we have \[
(\nabla_{\mathbf{X}}\boldsymbol{\Phi})(\mathbf{S})=\nabla_{\mathbf{X}}(\boldsymbol{\Phi}(\mathbf{S}))-\boldsymbol{\Phi}(\nabla_{\mathbf{X}} \mathbf{S}).
\] Applying this to \(\mathbf{S}=\mathbf{u}_1\otimes\cdots\otimes \mathbf{u}_r\) gives \[
\nabla_{\mathbf{X}} \mathbf{T}=(\nabla_{\mathbf{X}}\boldsymbol{\Phi})(\mathbf{u}_1\otimes\cdots\otimes \mathbf{u}_r)+\boldsymbol{\Phi}\bigl(\nabla_{\mathbf{X}}(\mathbf{u}_1\otimes\cdots\otimes \mathbf{u}_r)\bigr).
\]

By construction, \(\boldsymbol{\Phi}\) is a combination of tensor products, permutations, contractions, and isomorphisms induced by the metrics. Compatibility of the connections ensures that the musical identifications are parallel; contractions, permutations, and tensor products are also parallel for the induced connections. Each operation therefore commutes with the connection in the sense of Definition~\ref{def: homomorfismos paralelos}. It follows that \(\boldsymbol{\Phi}\) is parallel, that is, \(\nabla_{\mathbf{X}}\boldsymbol{\Phi}=0\).

Consequently, \[
\nabla_{\mathbf{X}} \mathbf{T}=\boldsymbol{\Phi}\bigl(\nabla_{\mathbf{X}}(\mathbf{u}_1\otimes\cdots\otimes \mathbf{u}_r)\bigr).
\] Applying the Leibniz rule for the tensor product, \[
\nabla_{\mathbf{X}}(\mathbf{u}_1\otimes\cdots\otimes \mathbf{u}_r)
=
\displaystyle\sum_{i=1}^{r}
\mathbf{u}_1\otimes\cdots\otimes (\nabla_{\mathbf{X}} \mathbf{u}_i)\otimes\cdots\otimes \mathbf{u}_r,
\] and substituting gives \[
\nabla_{\mathbf{X}} \mathbf{T}
=
\displaystyle\sum_{i=1}^{r}
\boldsymbol{\Phi}\bigl(\mathbf{u}_1\otimes\cdots\otimes (\nabla_{\mathbf{X}} \mathbf{u}_i)\otimes\cdots\otimes \mathbf{u}_r\bigr).
\] This is equivalent to \[
\nabla_{\mathbf{X}} \mathbf{T}
=
\displaystyle\sum_{i=1}^{r}
\mathbf{u}_1*\cdots *(\nabla_{\mathbf{X}} \mathbf{u}_i)*\cdots *\mathbf{u}_r,
\] as required. \end{proof} To study equivalence of norms, we need to control the commutator of covariant derivatives with the Laplacian. Recall the curvature $\mathbf{R}^{\mathbf{E}}$ of the connection, introduced in Definition~\ref{def:operadores-diferenciales-en-haces-curvatura-de-una-conexion}. The product connection on tensor bundles combines $\mathbf{R}^{\mathbf{E}}$ with the Riemannian curvature; both therefore appear in the following commutator. \begin{lemma}[Commutator of the Bochner Laplacian with covariant derivatives] \label{lema conmutador laplaciano orden m} \index{commutator of the Laplacian with covariant derivatives} Let $(M,\mathbf{g})$ be a Riemannian manifold with or without boundary, and let $\mathbf{E}\to M$ be a smooth vector bundle equipped with a bundle metric $\mathbf{h}_{\mathbf{E}}$ and a compatible connection $\nabla^{\mathbf{E}}$. For $\ell\in\mathbb N_0$, set \[
\mathbf{E}_\ell:=T^{(0,\ell)}(TM)\otimes\mathbf{E},
\] with the induced product metric and connection, and denote the covariant derivative by $\nabla_\ell\colon\Gamma(\mathbf{E}_\ell)\to\Gamma(\mathbf{E}_{\ell+1})$. For $m\in\mathbb N$ we write \[
\nabla^m:=\nabla_{\ell+m-1}\circ\cdots\circ\nabla_\ell,
\] so that each new derivative index occupies the first covariant factor. Let \[
\Delta_{B,\ell}:=(\nabla_\ell)_h^*\nabla_\ell=-\operatorname{tr}_{\mathbf{g}}(\nabla_\ell^2)
\] be the Bochner Laplacian on $\mathbf{E}_\ell$. Then, for every $\mathbf{u}\in\Gamma(\mathbf{E}_\ell)$, \begin{equation}
\label{eq:conmutador-bochner-orden-m-cap9}
\Delta_{B,\ell+m}(\nabla^m\mathbf{u})-\nabla^m(\Delta_{B,\ell}\mathbf{u})
=\displaystyle\sum_{k=0}^{m}\nabla^k\operatorname{Rm}*\nabla^{m-k}\mathbf{u}.
\end{equation} Here $\operatorname{Rm}$ represents either $\mathbf{R}^M$ or $\mathbf{R}^{\mathbf{E}}$, depending on the summand, with their induced connections, and $*$ denotes a finite sum of tensor products, permutations, and universal contractions. The coefficients and the number of terms depend only on $m$, $\ell$, $\dim M$, and the rank of $\mathbf{E}$. \end{lemma}

\begin{proof} We first prove the identity for a single derivative. Let $\mathbf{F}\to M$ be a smooth vector bundle equipped with a bundle metric $\mathbf{h}_{\mathbf{F}}$ and a compatible connection $\nabla^{\mathbf{F}}$, curvature $\mathbf{R}^{\mathbf{F}}$, and Bochner Laplacian \[
\Delta_B^{\mathbf{F}}:=(\nabla^{\mathbf{F}})^*\nabla^{\mathbf{F}}
=-\operatorname{tr}_{\mathbf{g}}\bigl((\nabla^{\mathbf{F}})^2\bigr).
\] We claim that, for each $\mathbf{v}\in\Gamma(\mathbf{F})$, each $x\in M$, each $\mathbf{X}\in T_xM$, and any orthonormal basis $(\mathbf{e}_1,\dots,\mathbf{e}_n)$ of $T_xM$, we have \begin{equation}
\label{eq:conmutador-bochner-primer-orden-exacto-cap9}
\begin{aligned}
\bigl(\nabla^{\mathbf{F}}\Delta_B^{\mathbf{F}}\mathbf{v}
-\Delta_B^{T^*M\otimes\mathbf{F}}(\nabla^{\mathbf{F}}\mathbf{v})\bigr)(\mathbf{X})=\nabla^{\mathbf{F}}_{\operatorname{Ric}^{\sharp}\mathbf{X}}\mathbf{v}
+2\displaystyle\sum_{i=1}^{n}\mathbf{R}^{\mathbf{F}}(\mathbf{e}_i,\mathbf{X})
\nabla^{\mathbf{F}}_{\mathbf{e}_i}\mathbf{v}
+\displaystyle\sum_{i=1}^{n}(\nabla_{\mathbf{e}_i}\mathbf{R}^{\mathbf{F}})
(\mathbf{e}_i,\mathbf{X})\mathbf{v}.
\end{aligned}
\end{equation} All terms on the right-hand side are evaluated at $x$. The sums are tensor contractions and thus do not depend on the chosen orthonormal basis.

First suppose that $x\in\operatorname{Int}(M)$. Both sides of \eqref{eq:conmutador-bochner-primer-orden-exacto-cap9} belong to $T_x^*M\otimes\mathbf{F}_x$; we therefore fix an arbitrary vector $\mathbf{X}\in T_xM$ and compare their values at $\mathbf{X}$. Extend $\mathbf{X},\mathbf{e}_1,\dots,\mathbf{e}_n$ to local fields, denoted in the same way, such that \[
\nabla^M_{\mathbf{Y}}\mathbf{X}(x)=0,
\qquad \nabla^M_{\mathbf{Y}}\mathbf{e}_i(x)=0
\] for every $\mathbf{Y}\in T_xM$ and every $i$. These extensions are obtained by taking constant components in normal coordinates centered at $x$.

Set $\mathbf{A}:=\nabla^{\mathbf{F}}\mathbf{v}\in
\Gamma(T^*M\otimes\mathbf{F})$. According to the convention in Proposition~\ref{prop:formulas-segunda-tercera-derivada-covariante}, the first two arguments of $\nabla^3\mathbf{v}=\nabla^2\mathbf{A}$ are the added derivative indices, and the third is the original argument of $\mathbf{A}$. Since the Levi--Civita connection commutes with metric contractions by Proposition~\ref{prop:conexion-conmuta-contracciones}, at the point $x$ we have \begin{equation}
\label{eq:dos-trazas-tercera-derivada-cap9}
\begin{aligned}
(\nabla^{\mathbf{F}}\Delta_B^{\mathbf{F}}\mathbf{v})(\mathbf{X})
&=-\displaystyle\sum_{i=1}^{n}(\nabla^3\mathbf{v})(\mathbf{X},\mathbf{e}_i,\mathbf{e}_i),\\
\bigl(\Delta_B^{T^*M\otimes\mathbf{F}}\mathbf{A}\bigr)(\mathbf{X})
&=-\displaystyle\sum_{i=1}^{n}(\nabla^3\mathbf{v})(\mathbf{e}_i,\mathbf{e}_i,\mathbf{X}).
\end{aligned}
\end{equation} The first equality follows by differentiating the tensor $-\operatorname{tr}_{\mathbf{g}}\nabla^2\mathbf{v}$, rather than a formula valid only at $x$. Subtracting the two equalities in \eqref{eq:dos-trazas-tercera-derivada-cap9} gives \begin{equation}
\label{eq:diferencia-trazas-tercera-derivada-cap9}
\begin{aligned}
\bigl(\nabla^{\mathbf{F}}\Delta_B^{\mathbf{F}}\mathbf{v}
-\Delta_B^{T^*M\otimes\mathbf{F}}\mathbf{A}\bigr)(\mathbf{X})=\displaystyle\sum_{i=1}^{n}\bigl[(\nabla^3\mathbf{v})(\mathbf{e}_i,\mathbf{e}_i,\mathbf{X})
-(\nabla^3\mathbf{v})(\mathbf{X},\mathbf{e}_i,\mathbf{e}_i)\bigr].
\end{aligned}
\end{equation}

We interchange the indices in two steps. The Ricci identity in Proposition~\ref{prop:curvatura-haz-tensorial-identidad-ricci}, first applied to $\mathbf{v}$ and then covariantly differentiated in the direction $\mathbf{e}_i$, gives \begin{equation}
\label{eq:primer-intercambio-tercera-derivada-cap9}
\begin{aligned}
(\nabla^3\mathbf{v})(\mathbf{e}_i,\mathbf{e}_i,\mathbf{X})
-(\nabla^3\mathbf{v})(\mathbf{e}_i,\mathbf{X},\mathbf{e}_i)=(\nabla_{\mathbf{e}_i}\mathbf{R}^{\mathbf{F}})(\mathbf{e}_i,\mathbf{X})\mathbf{v}
+\mathbf{R}^{\mathbf{F}}(\mathbf{e}_i,\mathbf{X})
\nabla^{\mathbf{F}}_{\mathbf{e}_i}\mathbf{v}.
\end{aligned}
\end{equation} Indeed, the terms in which the derivative acts on $\mathbf{e}_i$ or $\mathbf{X}$ vanish at $x$ by the choice of extensions.

We now apply the same Ricci identity to the section $\mathbf{A}=\nabla^{\mathbf{F}}\mathbf{v}$ of the bundle $T^*M\otimes\mathbf{F}$. By the product curvature formula \eqref{eq:curvatura-cotangente-producto-haz}, \begin{equation}
\label{eq:segundo-intercambio-tercera-derivada-cap9}
\begin{aligned}
(\nabla^3\mathbf{v})(\mathbf{e}_i,\mathbf{X},\mathbf{e}_i)
-(\nabla^3\mathbf{v})(\mathbf{X},\mathbf{e}_i,\mathbf{e}_i)=\mathbf{R}^{\mathbf{F}}(\mathbf{e}_i,\mathbf{X})
\nabla^{\mathbf{F}}_{\mathbf{e}_i}\mathbf{v}
-\nabla^{\mathbf{F}}_{\mathbf{R}^M(\mathbf{e}_i,\mathbf{X})\mathbf{e}_i}\mathbf{v}.
\end{aligned}
\end{equation}

Adding \eqref{eq:primer-intercambio-tercera-derivada-cap9} and \eqref{eq:segundo-intercambio-tercera-derivada-cap9}, then summing over $i$ and using \eqref{eq:diferencia-trazas-tercera-derivada-cap9}, gives \[
\begin{aligned}
\bigl(\nabla^{\mathbf{F}}\Delta_B^{\mathbf{F}}\mathbf{v}
-\Delta_B^{T^*M\otimes\mathbf{F}}(\nabla^{\mathbf{F}}\mathbf{v})\bigr)(\mathbf{X})=\displaystyle\sum_{i=1}^{n}(\nabla_{\mathbf{e}_i}\mathbf{R}^{\mathbf{F}})(\mathbf{e}_i,\mathbf{X})\mathbf{v}
+2\displaystyle\sum_{i=1}^{n}\mathbf{R}^{\mathbf{F}}(\mathbf{e}_i,\mathbf{X})
\nabla^{\mathbf{F}}_{\mathbf{e}_i}\mathbf{v}
-\nabla^{\mathbf{F}}_{\sum_{i=1}^{n}\mathbf{R}^M(\mathbf{e}_i,\mathbf{X})\mathbf{e}_i}\mathbf{v}.
\end{aligned}
\] Proposition~\ref{prop:contraccion-curvatura-endomorfismo-ricci} states that $-\displaystyle\sum_{i=1}^{n}\mathbf{R}^M(\mathbf{e}_i,\mathbf{X})\mathbf{e}_i
=\operatorname{Ric}^{\sharp}\mathbf{X}$, and we therefore obtain \eqref{eq:conmutador-bochner-primer-orden-exacto-cap9}. Since $x$ and $\mathbf{X}$ were arbitrary, the identity holds on $\operatorname{Int}(M)$. If $M$ has boundary, both sides are smooth sections up to $\partial M$ and $\operatorname{Int}(M)$ is dense in $M$; by continuity, the identity also holds on the boundary.

As a check on the sign, if $\mathbf{F}=M\times\mathbb R$ has the trivial connection, the preceding formula reduces to \[
d\Delta_B f-\Delta_B^{T^*M}(df)=df\circ\operatorname{Ric}^{\sharp}.
\]

We now apply the identity to $\mathbf{F}=\mathbf{E}_\ell=T^{(0,\ell)}(TM)\otimes\mathbf{E}$. By \eqref{eq:curvatura-haz-tensorial-general}, the curvature of $\mathbf{E}_\ell$ is the sum of the action of $\mathbf{R}^{\mathbf{E}}$ on the factor $\mathbf{E}$ and the actions of $\mathbf{R}^M$ on the covariant factors. Moving \eqref{eq:conmutador-bochner-primer-orden-exacto-cap9} to the other side and absorbing the signs and universal contractions into the symbol $*$ gives \begin{equation}
\label{eq:conmutador-bochner-primer-orden-cap9}
\Delta_{B,\ell+1}(\nabla\mathbf{u})-\nabla(\Delta_{B,\ell}\mathbf{u})
=\operatorname{Rm}*\nabla\mathbf{u}+\nabla\operatorname{Rm}*\mathbf{u}.
\end{equation} This is the case $m=1$.

Suppose that \eqref{eq:conmutador-bochner-orden-m-cap9} holds for some $m\geq1$. Then \begin{align*}
&\Delta_{B,\ell+m+1}(\nabla^{m+1}\mathbf{u})
-\nabla^{m+1}(\Delta_{B,\ell}\mathbf{u})\\
&\quad=\bigl(\Delta_{B,\ell+m+1}\nabla-\nabla\Delta_{B,\ell+m}\bigr)(\nabla^m\mathbf{u})\\
&\qquad+\nabla\bigl(\Delta_{B,\ell+m}(\nabla^m\mathbf{u})
-\nabla^m(\Delta_{B,\ell}\mathbf{u})\bigr).
\end{align*} By \eqref{eq:conmutador-bochner-primer-orden-cap9}, the first term is a universal sum of terms of the types $\operatorname{Rm}*\nabla^{m+1}\mathbf{u}$ and $\nabla\operatorname{Rm}*\nabla^m\mathbf{u}$. For the second we use the induction hypothesis and the Leibniz rule in Proposition~\ref{prop: leibniz operador *}: \begin{align*}
\nabla\bigl(\Delta_{B,\ell+m}(\nabla^m\mathbf{u})
-\nabla^m(\Delta_{B,\ell}\mathbf{u})\bigr)=\displaystyle\sum_{k=0}^{m}\Bigl(
\nabla^{k+1}\operatorname{Rm}*\nabla^{m-k}\mathbf{u}
+\nabla^k\operatorname{Rm}*\nabla^{m+1-k}\mathbf{u}\Bigr).
\end{align*} All the resulting terms have the form $\nabla^q\operatorname{Rm}*\nabla^{m+1-q}\mathbf{u}$ with $0\leq q\leq m+1$. Regrouping, for each $q$, the universal contractions that occur, we conclude that \[
\Delta_{B,\ell+m+1}(\nabla^{m+1}\mathbf{u})
-\nabla^{m+1}(\Delta_{B,\ell}\mathbf{u})
=\displaystyle\sum_{q=0}^{m+1}\nabla^q\operatorname{Rm}*\nabla^{m+1-q}\mathbf{u}.
\] This completes the induction. \end{proof}

\begin{theorem}[Equivalence of norms on $H^m(M,\mathbf{E})$ in terms of the Bochner Laplacian] \label{teo:sobolev_cerrada_riemanniana}\index{Sobolev space@Sobolev space!equivalent norm using the Bochner Laplacian} Let $(M,\mathbf{g})$ be a closed Riemannian manifold, and let $\mathbf{E}\to M$ be a smooth vector bundle equipped with a bundle metric $\mathbf{h}_{\mathbf{E}}$ and a compatible connection $\nabla^{\mathbf{E}}$. Let $m\in\mathbb N$. For $\mathbf{u}\in\Gamma(\mathbf{E})$ define \[
\|\mathbf{u}\|_{H^m(M,\mathbf{E})}^2
:=
\displaystyle\sum_{j=0}^{m}\int_M |\nabla^j \mathbf{u}|^2_{\mathbf{g},\mathbf{h}_{\mathbf{E}}}\,d\lambda_{\mathbf{g}},
\] and \[
\|\mathbf{u}\|_{\Delta_B}^2
:=
\displaystyle\sum_{\substack{j=0\\ j\text{ even}}}^{m}
\int_M |\Delta_B^{\frac{j}{2}}\mathbf{u}|^2_{\mathbf{h}_{\mathbf{E}}}\,d\lambda_{\mathbf{g}}
+
\displaystyle\sum_{\substack{j=1\\ j\text{ odd}}}^{m}
\int_M \big|\nabla \Delta_B^{\frac{j-1}{2}}\mathbf{u}\big|^2_{\mathbf{g},\mathbf{h}_{\mathbf{E}}}\,d\lambda_{\mathbf{g}}.
\] Denote by $H^m(M,\mathbf{E})$ the completion of $\Gamma(\mathbf{E})$ with respect to the first norm. Then there are constants $c,C>0$ such that \[
c\|\mathbf{u}\|_{H^m(M,\mathbf{E})}^2
\le
\|\mathbf{u}\|_{\Delta_B}^2
\le
C\|\mathbf{u}\|_{H^m(M,\mathbf{E})}^2,
\qquad \forall \mathbf{u}\in \Gamma(\mathbf{E}).
\] Consequently, the second norm determines the same completion. For an element of $H^m(M,\mathbf{E})$, each term containing a power of $\Delta_B$ is understood as the limit in $L^2(M,\mathbf E\oplus(T^*M\otimes\mathbf E))$ of the corresponding term of a smooth Cauchy sequence. \end{theorem}

\begin{proof} Since $M$ is closed, it is compact and $\partial M=\varnothing$. We prove the inequalities on $\Gamma(\mathbf{E})$; the assertion about completions will then follow by continuity.

Throughout what follows, all integrations by parts have no boundary terms.

\medskip

We first prove that \[
\|\mathbf{u}\|_{\Delta_B}^2 \le C \|\mathbf{u}\|_{H^m(M,\mathbf{E})}^2.
\]

By definition, \[
\Delta_B=\nabla^{*}\nabla=-\operatorname{tr}_{\mathbf{g}}(\nabla^2).
\] Since the connection is compatible with the metric, covariant differentiation commutes with contraction. Define the contractions \begin{align*}
C_k^{\mathrm{even}}(T)
&:=g^{i_1i_2}\cdots g^{i_{2k-1}i_{2k}}
T_{i_1\cdots i_{2k}},\\
\bigl(C_k^{\mathrm{odd}}(S)\bigr)_a
&:=g^{i_1i_2}\cdots g^{i_{2k-1}i_{2k}}
S_{a i_1\cdots i_{2k}}.
\end{align*} The first contracts the pairs $(1,2),(3,4),\ldots$; the second leaves the first index free and contracts $(2,3),(4,5),\ldots$. Induction on $k$ gives \[
\Delta_B^k \mathbf{u}=(-1)^kC_k^{\mathrm{even}}(\nabla^{2k}\mathbf{u}).
\] Set $\Phi_k:=(-1)^kC_k^{\mathrm{even}}$. Since $\Phi_k$ is obtained using contractions and linear combinations, Lemma~\ref{lema estimacion estrella} implies that there is a constant $C_k>0$ such that \[
|\Delta_B^k \mathbf{u}|_{\mathbf{h}_{\mathbf{E}}}
\le
C_k |\nabla^{2k}\mathbf{u}|_{\mathbf{g},\mathbf{h}_{\mathbf{E}}}.
\] Squaring and integrating gives \[
\int_M |\Delta_B^k \mathbf{u}|_{\mathbf{h}_{\mathbf{E}}}^2\,d\lambda_{\mathbf{g}}
\le
C_k^2 \int_M |\nabla^{2k}\mathbf{u}|_{\mathbf{g},\mathbf{h}_{\mathbf{E}}}^2\,d\lambda_{\mathbf{g}}.
\]

Since $\nabla \mathbf{g}=0$, differentiating each factor $\mathbf{g}^{-1}$ used in the contraction gives zero; by the Leibniz rule, \[
\nabla \Delta_B^k \mathbf{u}
=
(-1)^k \nabla\bigl(C_k^{\mathrm{even}}(\nabla^{2k}\mathbf{u})\bigr)
=
(-1)^k C_k^{\mathrm{odd}}(\nabla^{2k+1}\mathbf{u}).
\] We therefore define \[
\Psi_k\colon T^{(0,2k+1)}(TM)\otimes \mathbf{E}\longrightarrow T^*M\otimes \mathbf{E},
\qquad
\Psi_k(S):=(-1)^kC_k^{\mathrm{odd}}(S).
\] Thus \[
\nabla \Delta_B^k \mathbf{u}=\Psi_k(\nabla^{2k+1}\mathbf{u}),
\] and, again by Lemma~\ref{lema estimacion estrella}, there is $C_k'>0$ such that \[
|\nabla \Delta_B^k \mathbf{u}|_{\mathbf{g},\mathbf{h}_{\mathbf{E}}}
\le
C_k' |\nabla^{2k+1}\mathbf{u}|_{\mathbf{g},\mathbf{h}_{\mathbf{E}}}.
\] Squaring and integrating gives \[
\int_M |\nabla \Delta_B^k \mathbf{u}|_{\mathbf{g},\mathbf{h}_{\mathbf{E}}}^2\,d\lambda_{\mathbf{g}}
\le
(C_k')^2 \int_M |\nabla^{2k+1}\mathbf{u}|_{\mathbf{g},\mathbf{h}_{\mathbf{E}}}^2\,d\lambda_{\mathbf{g}}.
\]

Consequently, each even summand of $\|\mathbf{u}\|_{\Delta_B}^2$ is bounded by the corresponding even-order summand of $\|\mathbf{u}\|_{H^m(M,\mathbf{E})}^2$, and each odd summand of $\|\mathbf{u}\|_{\Delta_B}^2$ is bounded by the corresponding odd-order summand of $\|\mathbf{u}\|_{H^m(M,\mathbf{E})}^2$. Adding all the terms, we conclude that there is a constant $C>0$ such that \[
\|\mathbf{u}\|_{\Delta_B}^2 \le C \|\mathbf{u}\|_{H^m(M,\mathbf{E})}^2.
\]

\medskip

We now prove that \[
\|\mathbf{u}\|_{H^m(M,\mathbf{E})}^2 \le C \|\mathbf{u}\|_{\Delta_B}^2.
\]

We formulate the inductive estimate for all smooth sections. For $0\leq j\leq m$ and $\mathbf{v}\in\Gamma(\mathbf{E})$, introduce the truncated quantity \begin{equation}
\label{eq:norma-laplaciana-truncada-cap9}
\mathcal N_j(\mathbf{v})^2
:=
\displaystyle\sum_{\substack{0\leq r\leq j\\r\text{ even}}}
\|\Delta_B^{\frac{r}{2}}\mathbf{v}\|_{L^2(M,\mathbf{E})}^2
+
\displaystyle\sum_{\substack{1\leq r\leq j\\r\text{ odd}}}
\|\nabla^{\mathbf{E}}\Delta_B^{\frac{r-1}{2}}\mathbf{v}\|_{L^2(M,T^*M\otimes \mathbf{E})}^2.
\end{equation} Thus $\mathcal N_m(\mathbf{u})=\|\mathbf{u}\|_{\Delta_B}$. We prove the following quantified assertion by induction: for each $j\in\{0,\dots,m\}$ there is $C_j>0$ such that, for every smooth section $\mathbf{v}\in\Gamma(\mathbf{E})$, \[
\int_M |\nabla^j \mathbf{v}|_{\mathbf{g},\mathbf{h}_{\mathbf{E}}}^2\,d\lambda_{\mathbf{g}}
\le
C_j \mathcal N_j(\mathbf{v})^2.
\] Since the estimate is formulated for every smooth section, at the step of order $j-2$ we may take $\mathbf{v}=\Delta_B\mathbf{u}$.

For $j=0$ and every $\mathbf{v}\in\Gamma(\mathbf{E})$ we have \[
\int_M |\mathbf{v}|_{\mathbf{h}_{\mathbf{E}}}^2\,d\lambda_{\mathbf{g}}
=
\mathcal N_0(\mathbf{v})^2,
\] since the truncated sum contains only the term of order zero.

For $j=1$ and every $\mathbf{v}\in\Gamma(\mathbf{E})$ we have \[
\int_M |\nabla \mathbf{v}|_{\mathbf{g},\mathbf{h}_{\mathbf{E}}}^2\,d\lambda_{\mathbf{g}}
\leq
\mathcal N_1(\mathbf{v})^2,
\] since this is the odd term of order one.

Now consider the case $j=2$. Applying Corollary~\ref{cor:green-bochner} to the bundle $T^*M\otimes \mathbf{E}$ gives \[
\int_M |\nabla^2\mathbf{u}|_{\mathbf{g},\mathbf{h}_{\mathbf{E}}}^2\,d\lambda_{\mathbf{g}}
=
\int_M \langle \nabla \mathbf{u},\Delta_{B,1}(\nabla \mathbf{u})\rangle_{\mathbf{g},\mathbf{h}_{\mathbf{E}}}\,d\lambda_{\mathbf{g}}.
\] By Lemma~\ref{lema conmutador laplaciano orden m} with $m=1$, \[
\Delta_{B,1}(\nabla \mathbf{u})
=
\nabla(\Delta_{B,0} \mathbf{u})
+\operatorname{Rm}*\nabla \mathbf{u}
+\nabla\operatorname{Rm}*\mathbf{u}.
\] Substituting, \[
\begin{aligned}
\int_M |\nabla^2\mathbf{u}|_{\mathbf{g},\mathbf{h}_{\mathbf{E}}}^2\,d\lambda_{\mathbf{g}}
={}&
\int_M \langle \nabla \mathbf{u},\nabla(\Delta_B \mathbf{u})\rangle_{\mathbf{g},\mathbf{h}_{\mathbf{E}}}\,d\lambda_{\mathbf{g}}
+\int_M \langle \nabla \mathbf{u},\operatorname{Rm}*\nabla \mathbf{u}\rangle_{\mathbf{g},\mathbf{h}_{\mathbf{E}}}\,d\lambda_{\mathbf{g}}\\
&+\int_M \langle \nabla \mathbf{u},\nabla\operatorname{Rm}*\mathbf{u}\rangle_{\mathbf{g},\mathbf{h}_{\mathbf{E}}}\,d\lambda_{\mathbf{g}}.
\end{aligned}
\] Applying Corollary~\ref{cor:green-bochner} once more, now to the bundle $\mathbf{E}$, gives \[
\int_M \langle \nabla \mathbf{u},\nabla(\Delta_B \mathbf{u})\rangle_{\mathbf{g},\mathbf{h}_{\mathbf{E}}}\,d\lambda_{\mathbf{g}}
=
\int_M |\Delta_B \mathbf{u}|_{\mathbf{h}_{\mathbf{E}}}^2\,d\lambda_{\mathbf{g}}.
\] Moreover, since $M$ is compact, the curvature tensor is uniformly bounded and, by Lemma~\ref{lema estimacion estrella}, there is a constant $C>0$ such that \[
|\operatorname{Rm}*\nabla \mathbf{u}|_{\mathbf{g},\mathbf{h}_{\mathbf{E}}}
\le
C |\nabla \mathbf{u}|_{\mathbf{g},\mathbf{h}_{\mathbf{E}}}.
\] Consequently, \[
\left|
\int_M \langle \nabla \mathbf{u},\operatorname{Rm}*\nabla \mathbf{u}\rangle_{\mathbf{g},\mathbf{h}_{\mathbf{E}}}\,d\lambda_{\mathbf{g}}
\right|
\leq
C\int_M |\nabla \mathbf{u}|_{\mathbf{g},\mathbf{h}_{\mathbf{E}}}^2\,d\lambda_{\mathbf{g}}.
\] Likewise, the compactness of $M$ gives a uniform bound on $\nabla\operatorname{Rm}$, so that \[
|\nabla\operatorname{Rm}*\mathbf{u}|_{\mathbf{g},\mathbf{h}_{\mathbf{E}}}\leq C|\mathbf{u}|_{\mathbf{h}_{\mathbf{E}}}.
\] By the Cauchy--Schwarz inequality, \[
\left|\int_M\langle \nabla \mathbf{u},\nabla\operatorname{Rm}*\mathbf{u}\rangle_{\mathbf{g},\mathbf{h}_{\mathbf{E}}}\,d\lambda_{\mathbf{g}}\right|
\leq C\|\nabla \mathbf{u}\|_{L^2(M,T^*M\otimes \mathbf{E})}\|\mathbf{u}\|_{L^2(M,\mathbf{E})}.
\] The two factors correspond to lower-order terms appearing in $\|\mathbf{u}\|_{\Delta_B}$. The elementary inequality $2ab\leq a^2+b^2$, applied with $a=\|\nabla \mathbf{u}\|_{L^2(M,T^*M\otimes \mathbf{E})}$ and $b=\|\mathbf{u}\|_{L^2(M,\mathbf{E})}$, gives \[
\left|\int_M\langle \nabla \mathbf{u},\nabla\operatorname{Rm}*\mathbf{u}\rangle_{\mathbf{g},\mathbf{h}_{\mathbf{E}}}\,d\lambda_{\mathbf{g}}\right|
\leq C\|\nabla \mathbf{u}\|_{L^2(M,T^*M\otimes \mathbf{E})}^2+C\|\mathbf{u}\|_{L^2(M,\mathbf{E})}^2.
\] Consequently, \[
\int_M |\nabla^2\mathbf{u}|_{\mathbf{g},\mathbf{h}_{\mathbf{E}}}^2\,d\lambda_{\mathbf{g}}
\leq
\int_M |\Delta_B \mathbf{u}|_{\mathbf{h}_{\mathbf{E}}}^2\,d\lambda_{\mathbf{g}}
+C\int_M |\nabla \mathbf{u}|_{\mathbf{g},\mathbf{h}_{\mathbf{E}}}^2\,d\lambda_{\mathbf{g}}
+C\int_M |\mathbf{u}|_{\mathbf{h}_{\mathbf{E}}}^2\,d\lambda_{\mathbf{g}}.
\] The three terms on the right-hand side appear, respectively, in $\mathcal N_2(\mathbf{u})^2$. Therefore, \[
\int_M |\nabla^2\mathbf{u}|_{\mathbf{g},\mathbf{h}_{\mathbf{E}}}^2\,d\lambda_{\mathbf{g}}
\le
C_2\mathcal N_2(\mathbf{u})^2.
\]

Now let $j\geq 3$ and suppose that, for every integer $l<j$ and every $\mathbf{v}\in\Gamma(\mathbf{E})$, we have already proved that \[
\int_M |\nabla^l \mathbf{v}|_{\mathbf{g},\mathbf{h}_{\mathbf{E}}}^2\,d\lambda_{\mathbf{g}}
\le
C_l\mathcal N_l(\mathbf{v})^2.
\] We want to prove the same estimate for $j$.

Write $\nabla^j\mathbf{u}=\nabla^2(\nabla^{j-2}\mathbf{u})$. The argument for the case $j=2$ depends only on compatibility of the connection and the uniform bounds on the curvature and its first derivative. It therefore applies to the bundle $T^{(0,j-2)}(TM)\otimes \mathbf{E}$ and the section $\nabla^{j-2}\mathbf{u}$. The curvature of this product bundle and its first derivative are sums of contractions of $\mathbf{R}^M$, $\mathbf{R}^{\mathbf{E}}$, and their derivatives; all are bounded because $M$ is compact. This yields \[
\begin{aligned}
\int_M |\nabla^j \mathbf{u}|_{\mathbf{g},\mathbf{h}_{\mathbf{E}}}^2\,d\lambda_{\mathbf{g}}
\leq{}&
\int_M |\Delta_{B,j-2}(\nabla^{j-2}\mathbf{u})|_{\mathbf{g},\mathbf{h}_{\mathbf{E}}}^2\,d\lambda_{\mathbf{g}}
+C\int_M |\nabla^{j-1}\mathbf{u}|_{\mathbf{g},\mathbf{h}_{\mathbf{E}}}^2\,d\lambda_{\mathbf{g}}\\
&+C\int_M |\nabla^{j-2}\mathbf{u}|_{\mathbf{g},\mathbf{h}_{\mathbf{E}}}^2\,d\lambda_{\mathbf{g}}.
\end{aligned}
\] By the induction hypothesis, the last two terms are bounded by $C\mathcal N_j(\mathbf{u})^2$, since $\mathcal N_l(\mathbf{u})\leq\mathcal N_j(\mathbf{u})$ for $l<j$. Thus it suffices to estimate the first.

Applying Lemma~\ref{lema conmutador laplaciano orden m} with $m=j-2$ gives \[
\Delta_{B,j-2}(\nabla^{j-2}\mathbf{u})
=
\nabla^{j-2}(\Delta_B \mathbf{u})
+
\displaystyle\sum_{k=0}^{j-2}\nabla^k\operatorname{Rm}*\nabla^{j-2-k}\mathbf{u}.
\] Squaring, using the inequality \[
|A_0+\cdots+A_N|^2\le C\displaystyle\sum_{r=0}^N |A_r|^2,
\], and integrating, we obtain \[
\int_M |\Delta_{B,j-2}(\nabla^{j-2}\mathbf{u})|_{\mathbf{g},\mathbf{h}_{\mathbf{E}}}^2\,d\lambda_{\mathbf{g}}
\le
C\int_M |\nabla^{j-2}(\Delta_B \mathbf{u})|_{\mathbf{g},\mathbf{h}_{\mathbf{E}}}^2\,d\lambda_{\mathbf{g}}
+
C\displaystyle\sum_{k=0}^{j-2}\int_M |\nabla^k\operatorname{Rm}*\nabla^{j-2-k}\mathbf{u}|_{\mathbf{g},\mathbf{h}_{\mathbf{E}}}^2\,d\lambda_{\mathbf{g}}.
\] Since $M$ is compact, all covariant derivatives of \(\operatorname{Rm}\) are uniformly bounded. By Lemma~\ref{lema estimacion estrella}, for each \(k\) there is a constant \(C_k>0\) such that \[
|\nabla^k\operatorname{Rm}*\nabla^{j-2-k}\mathbf{u}|_{\mathbf{g},\mathbf{h}_{\mathbf{E}}}
\le
C_k |\nabla^{j-2-k}\mathbf{u}|_{\mathbf{g},\mathbf{h}_{\mathbf{E}}}.
\] Therefore, \[
\int_M |\Delta_{B,j-2}(\nabla^{j-2}\mathbf{u})|_{\mathbf{g},\mathbf{h}_{\mathbf{E}}}^2\,d\lambda_{\mathbf{g}}
\le
C\int_M |\nabla^{j-2}(\Delta_B \mathbf{u})|_{\mathbf{g},\mathbf{h}_{\mathbf{E}}}^2\,d\lambda_{\mathbf{g}}
+
C\displaystyle\sum_{l=0}^{j-2}\int_M |\nabla^l \mathbf{u}|_{\mathbf{g},\mathbf{h}_{\mathbf{E}}}^2\,d\lambda_{\mathbf{g}}.
\] By the induction hypothesis, the lower-order sum is controlled by $\mathcal N_j(\mathbf{u})^2$. Consequently, \[
\int_M |\nabla^j \mathbf{u}|_{\mathbf{g},\mathbf{h}_{\mathbf{E}}}^2\,d\lambda_{\mathbf{g}}
\le
C\int_M |\nabla^{j-2}(\Delta_B \mathbf{u})|_{\mathbf{g},\mathbf{h}_{\mathbf{E}}}^2\,d\lambda_{\mathbf{g}}
+
C\mathcal N_j(\mathbf{u})^2.
\]

We now apply the inductive assertion, which holds for every smooth section, to $\mathbf{v}=\Delta_B\mathbf{u}$ at order $j-2$. Since $j-2<j$, we obtain \[
\int_M |\nabla^{j-2}(\Delta_B \mathbf{u})|_{\mathbf{g},\mathbf{h}_{\mathbf{E}}}^2\,d\lambda_{\mathbf{g}}
\le
C\left(
\displaystyle\sum_{\substack{\ell=0\\ \ell\text{ even}}}^{j-2}
\int_M \bigl|\Delta_B^{\frac{\ell}{2}}(\Delta_B \mathbf{u})\bigr|_{\mathbf{h}_{\mathbf{E}}}^2\,d\lambda_{\mathbf{g}}
+
\displaystyle\sum_{\substack{\ell=1\\ \ell\text{ odd}}}^{j-2}
\int_M \bigl|\nabla \Delta_B^{\frac{\ell-1}{2}}(\Delta_B \mathbf{u})\bigr|_{\mathbf{g},\mathbf{h}_{\mathbf{E}}}^2\,d\lambda_{\mathbf{g}}
\right).
\] We now rewrite each of these sums.

For the even sum, \[
\Delta_B^{\frac{\ell}{2}}(\Delta_B \mathbf{u})=\Delta_B^{\frac{\ell}{2}+1}\mathbf{u}.
\] If \(\ell\) ranges over the even integers between \(0\) and \(j-2\), then \(r=\ell+2\) ranges over the even integers between \(2\) and \(j\). Therefore, \[
\displaystyle\sum_{\substack{\ell=0\\ \ell\text{ even}}}^{j-2}
\int_M \bigl|\Delta_B^{\frac{\ell}{2}}(\Delta_B \mathbf{u})\bigr|_{\mathbf{h}_{\mathbf{E}}}^2\,d\lambda_{\mathbf{g}}
=
\displaystyle\sum_{\substack{r=2\\ r\text{ even}}}^{j}
\int_M \bigl|\Delta_B^{\frac{r}{2}}\mathbf{u}\bigr|_{\mathbf{h}_{\mathbf{E}}}^2\,d\lambda_{\mathbf{g}}.
\]

For the odd sum, \[
\nabla \Delta_B^{\frac{\ell-1}{2}}(\Delta_B \mathbf{u})
=
\nabla \Delta_B^{\frac{\ell+1}{2}}\mathbf{u}.
\] If \(\ell\) ranges over the odd integers between \(1\) and \(j-2\), then \(s=\ell+2\) ranges over the odd integers between \(3\) and \(j\). Thus \[
\displaystyle\sum_{\substack{\ell=1\\ \ell\text{ odd}}}^{j-2}
\int_M \bigl|\nabla \Delta_B^{\frac{\ell-1}{2}}(\Delta_B \mathbf{u})\bigr|_{\mathbf{g},\mathbf{h}_{\mathbf{E}}}^2\,d\lambda_{\mathbf{g}}
=
\displaystyle\sum_{\substack{s=3\\ s\text{ odd}}}^{j}
\int_M \bigl|\nabla \Delta_B^{\frac{s-1}{2}}\mathbf{u}\bigr|_{\mathbf{g},\mathbf{h}_{\mathbf{E}}}^2\,d\lambda_{\mathbf{g}}.
\]

Substituting these identities into the preceding estimate gives \[
\int_M |\nabla^{j-2}(\Delta_B \mathbf{u})|_{\mathbf{g},\mathbf{h}_{\mathbf{E}}}^2\,d\lambda_{\mathbf{g}}
\le
C\left(
\displaystyle\sum_{\substack{r=2\\ r\text{ even}}}^{j}
\int_M \bigl|\Delta_B^{\frac{r}{2}}\mathbf{u}\bigr|_{\mathbf{h}_{\mathbf{E}}}^2\,d\lambda_{\mathbf{g}}
+
\displaystyle\sum_{\substack{s=3\\ s\text{ odd}}}^{j}
\int_M \bigl|\nabla \Delta_B^{\frac{s-1}{2}}\mathbf{u}\bigr|_{\mathbf{g},\mathbf{h}_{\mathbf{E}}}^2\,d\lambda_{\mathbf{g}}
\right).
\] Each of these terms appears in the definition of $\mathcal N_j(\mathbf{u})^2$. Equivalently, $\mathcal N_{j-2}(\Delta_B\mathbf{u})\leq\mathcal N_j(\mathbf{u})$: the even powers and odd terms shift by exactly two orders. Consequently, \[
\int_M |\nabla^{j-2}(\Delta_B \mathbf{u})|_{\mathbf{g},\mathbf{h}_{\mathbf{E}}}^2\,d\lambda_{\mathbf{g}}
\le
C\mathcal N_j(\mathbf{u})^2.
\] Returning to the previous estimate for \(\displaystyle\int_M |\nabla^j \mathbf{u}|_{\mathbf{g},\mathbf{h}_{\mathbf{E}}}^2\,d\lambda_{\mathbf{g}}\), we conclude that \[
\int_M |\nabla^j \mathbf{u}|_{\mathbf{g},\mathbf{h}_{\mathbf{E}}}^2\,d\lambda_{\mathbf{g}}
\le
C_j\mathcal N_j(\mathbf{u})^2.
\] This completes the induction step. Finally, $\mathcal N_j(\mathbf{u})\leq\mathcal N_m(\mathbf{u})=\|\mathbf{u}\|_{\Delta_B}$ for $j\leq m$. Summing over $j\in\{0,\dots,m\}$ gives \[
\|\mathbf{u}\|_{H^m(M,\mathbf{E})}^2
=
\displaystyle\sum_{j=0}^{m}\int_M |\nabla^j \mathbf{u}|_{\mathbf{g},\mathbf{h}_{\mathbf{E}}}^2\,d\lambda_{\mathbf{g}}
\le
C\|\mathbf{u}\|_{\Delta_B}^2.
\]

Together with the inequality proved at the beginning, \[
\|\mathbf{u}\|_{\Delta_B}^2 \le C\|\mathbf{u}\|_{H^m(M,\mathbf{E})}^2,
\] this gives the two stated inequalities with positive constants. If $(\mathbf{u}_\nu)$ is Cauchy in the norm of $H^m(M,\mathbf E)$, the forward inequality shows that each sequence $\Delta_B^k\mathbf{u}_\nu$ and $\nabla\Delta_B^k\mathbf{u}_\nu$ appearing in the structured norm is Cauchy in $L^2(M,\mathbf E)$ and in $L^2(M,T^*M\otimes\mathbf E)$, respectively. We define these terms in the completion by their limits. The reverse inequality proves that this procedure adds no new elements and is independent of the representative sequence. This justifies the final assertion in the statement. \end{proof}

\subsection{Local characterization of Sobolev spaces in vector bundles} In this subsection we show how to characterize Sobolev sections through their component functions. To this end, we first prove an analogue of Lemma~\ref{lema cg} for arbitrary vector bundles: \begin{lemma}\label{lema equivalencia metricas fibradas} Let $M$ be a smooth manifold with or without boundary, let $K\subseteq M$ be compact, and let $\mathbf{E}\longrightarrow M$ be a smooth vector bundle of finite rank. Let $\mathbf{h}$ and $\widetilde{\mathbf{h}}$ be two smooth bundle metrics on $\mathbf{E}$. Then there are constants $c,C>0$ such that, for every $p\in K$ and every $v\in \mathbf{E}_{p}$, \[
c\,\mathbf{h}_{p}(v,v)\leq \widetilde{\mathbf{h}}_{p}(v,v)\leq C\,\mathbf{h}_{p}(v,v).
\] Equivalently, $\sqrt c\,|v|_{\mathbf h}\leq|v|_{\widetilde{\mathbf h}}
\leq\sqrt C\,|v|_{\mathbf h}$ in the same fibers. In particular, if $M$ is compact, we may take $K=M$. \end{lemma}

\begin{proof} If $K=\varnothing$ or the bundle has rank zero, it suffices to take $c=C=1$. Suppose that $K\neq\varnothing$ and the rank is $r\geq1$. Let $p\in M$. Consider \[
S_{p}:=\{v\in \mathbf{E}_{p}\mid \mathbf{h}_{p}(v,v)=1\}.
\]

This set is closed and bounded in a finite-dimensional normed vector space, and is therefore compact. Moreover, the function $f_{p}\colon S_{p}\longrightarrow\mathbb{R}$, given by $f_{p}(v)=\widetilde{\mathbf{h}}_{p}(v,v)$, is continuous and hence attains its maximum. Thus the function $f\colon M\longrightarrow\mathbb{R}$ defined by \[
f(p):=\max_{v\in S_{p}}f_{p}(v)
\] is well defined. We prove its continuity by identifying the spheres in nearby fibers with a fixed sphere. On a neighborhood $U$ of any point $p_*\in M$, the Gram--Schmidt process applied to a smooth frame gives a frame $\mathbf e_1,\ldots,\mathbf e_r$ orthonormal for $\mathbf h$. If the bundle is real, set $\mathbb K=\mathbb R$; if it is complex with Hermitian metrics, set $\mathbb K=\mathbb C$. Define \[
 S:=\left\{w\in\mathbb K^r\ \middle|\ \sum_{a=1}^r|w^a|^2=1\right\},
 \qquad
 F(q,w):=\widetilde{\mathbf h}_q
 \left(\sum_{a=1}^{r} w^a\mathbf e_a(q),\sum_{a=1}^{r} w^a\mathbf e_a(q)\right).
\] For $q\in U$, the map $w\mapsto\displaystyle\sum_{a=1}^{r} w^a\mathbf e_a(q)$ identifies $S$ with $S_q$. Hence $f(q)=\displaystyle\max_{w\in S}F(q,w)$. Let $H_{ab}(q)=\widetilde{\mathbf h}_q(\mathbf e_a(q),\mathbf e_b(q))$ be the continuous coefficients of the second metric in that frame. The Cauchy--Schwarz inequality gives $\displaystyle\sum_{a=1}^{r}|w^a|\leq\sqrt r$ for $w\in S$. Consequently, \[
 \begin{aligned}
 \sup_{w\in S}|F(q,w)-F(p_*,w)|
 &\leq \max_{1\leq a,b\leq r}|H_{ab}(q)-H_{ab}(p_*)|
       \sup_{w\in S}\left(\sum_{a=1}^{r}|w^a|\right)^2\\
 &\leq r\max_{1\leq a,b\leq r}|H_{ab}(q)-H_{ab}(p_*)|
 \longrightarrow0
 \qquad(q\to p_*).
 \end{aligned}
\] If two functions on $S$ differ by at most $\varepsilon$ at every point, their maxima differ by at most $\varepsilon$: first use $F(q,w)\leq F(p_*,w)+\varepsilon$, and then the inequality with the two points interchanged. Thus \[
 |f(q)-f(p_*)|
 \leq\sup_{w\in S}|F(q,w)-F(p_*,w)|,
\] which proves continuity of $f$ on $p_*$. Since the point was arbitrary, $f$ is continuous on $M$. The compactness of $K$ gives $p_0\in K$ such that $f(p)\leq f(p_0)=:C$ for every $p\in K$.

Moreover, \[
C=f(p_{0})=\max_{v\in S_{p_{0}}}\widetilde{\mathbf{h}}_{p_{0}}(v,v)>0,
\] since $v\neq 0$ for every $v\in S_{p_{0}}$. Thus, for every $p\in K$ and every $v\in S_{p}$, we have \[
\widetilde{\mathbf{h}}_{p}(v,v)\leq C.
\]

Define $\displaystyle m(p):=\min_{v\in S_p}\widetilde{\mathbf{h}}_p(v,v)$. In the preceding frame, $m(q)=\min_{w\in S}F(q,w)$. The same inequalities between $F(q,w)$ and $F(p_*,w)$, now taking minima, give \[
 |m(q)-m(p_*)|
 \leq\sup_{w\in S}|F(q,w)-F(p_*,w)|\longrightarrow0.
\] Thus $m$ is also continuous. Each $m(p)$ is positive because the minimum is attained at a nonzero vector and $\widetilde{\mathbf h}_p$ is positive definite. The compactness of $K$ implies that $m$ attains a minimum there, which remains positive: $\displaystyle c:=\min_{p\in K}m(p)>0$. Therefore, for every $p\in K$ and every $v\in S_p$, we have \[
c\leq \widetilde{\mathbf{h}}_{p}(v,v).
\]

Now let $p\in K$ and $v\in \mathbf{E}_{p}\setminus\{0\}$. Then $\displaystyle\frac{v}{|v|_{\mathbf{h}_{p}}}\in S_{p}$, so that \[
c\leq
\widetilde{\mathbf{h}}_{p}\!\left(\frac{v}{|v|_{\mathbf{h}_{p}}},\frac{v}{|v|_{\mathbf{h}_{p}}}\right)
\leq C.
\] This implies that \[
c\,|v|^{2}_{\mathbf{h}_{p}}
\leq \widetilde{\mathbf{h}}_{p}(v,v)
\leq C\,|v|^{2}_{\mathbf{h}_{p}},
\] that is, \[
c\,\mathbf{h}_{p}(v,v)\leq \widetilde{\mathbf{h}}_{p}(v,v)\leq C\,\mathbf{h}_{p}(v,v).
\] The case $v=0$ is immediate, since all the preceding quantities equal $0$. \end{proof}

With this in mind, we are ready to prove the following lemma, which allows us to work locally with these globally and intrinsically defined Sobolev spaces. Unlike the abstract definition of $W^{m,p}(M,\mathbf E)$, this lemma assumes that the connection is compatible with the bundle metric, because the proof uses the explicit form of the formal adjoint obtained from the integration-by-parts formula. \begin{lemma}\label{lema:meyers-serrin-haz-E} Let $(M,\mathbf{g})$ be a Riemannian manifold without boundary. Let $m\in\mathbb N_0$ and $1\leq p<\infty$, let $(U,\phi)$ be a regular coordinate ball in $M$, and let $\pi_{\mathbf{E}}\colon \mathbf{E}\longrightarrow M$ be a smooth real vector bundle of rank $r$ equipped with a bundle metric $\mathbf{h}_{\mathbf{E}}$ and a compatible connection $\nabla^{\mathbf{E}}$. Let $\mathbf{F}\in L^{1}_{\operatorname{loc}}(M,\mathbf{E})$ and let $\{\mathbf{e}_{1},\dots,\mathbf{e}_{r}\}$ be a local frame of $\mathbf{E}$ defined on a neighborhood of $\overline U$. Suppose that $\mathbf{F}=\displaystyle\sum_{a=1}^{r}F^{a}\mathbf{e}_{a}$ on $U$. Then \[
\mathbf{F}\restriction_{U}\in W^{m,p}(U,\mathbf{E}\restriction_{U}) \quad \text{if and only if} \quad
F^{a}\circ \phi^{-1}\in W^{m,p}(\phi(U))
\] for every $1\le a\le r$.

Moreover, there are constants $C_{m},c_{m}>0$ (independent of $\mathbf{F}$) such that \[c_{m}\displaystyle\sum_{a=1}^{r}
\left\|F^{a}\circ \phi^{-1}\right\|_{W^{m,p}(\phi(U))}\leq \|\mathbf{F}\|_{W^{m,p}(U,\mathbf{E}\restriction_{U})}
\leq
C_{m}\displaystyle\sum_{a=1}^{r}
\left\|F^{a}\circ \phi^{-1}\right\|_{W^{m,p}(\phi(U))}
\] for every $\mathbf{F}\in W^{m,p}(U,\mathbf{E}\restriction_{U})$. \end{lemma}

\begin{proof} The case $m=0$ follows directly from uniform equivalence of the fiber norms over the compact set $\overline U$ and equivalence of the volume elements. Thus suppose that $m\geq1$. Fix $s\ge 1$ and denote by $I_s=(i_{1},\dots,i_{s})$, $J_s=(j_{1},\dots,j_{s})$, etc., the ordered tuples or blocks of indices in $\{1,\dots,n\}$. We adopt the Einstein summation convention for these blocks. Define \[
\partial_{I_s} := \partial_{i_{1}}\cdots\partial_{i_{s}}, \quad
g^{I_s J_s} := \prod_{l=1}^{s}g^{i_{l}j_{l}}, \quad
g_{J_s T_s} := \prod_{l=1}^{s}g_{j_{l}t_{l}}.
\] Also, let $\beta_{I_s}\in \mathbb{N}_{0}^{n}$ be the multi-index associated with the block $I_s$ such that $\partial^{\beta_{I_s}} \equiv \partial_{I_s}$.

Since $\nabla^{\mathbf E}$ is compatible with $\mathbf h_{\mathbf E}$, we may use the expression for the adjoint obtained in Corollary~\ref{cor: expresion coordenada adjunto formal s derivada covariante}. Set $h_{ab}:=\mathbf h_{\mathbf E}(\mathbf e_a,\mathbf e_b)$ and $(h^{ab}):=(h_{ab})^{-1}$. Written in the chosen frame, with the leading term separated from the lower-order terms, that formula gives \[
\bigl((\nabla^s)^*\boldsymbol\Psi\bigr)^a
=(-1)^s\frac{h^{ab}}{\sqrt{\det(\mathbf{g})}}\partial_{I_s}\!\left(\sqrt{\det(\mathbf{g})} h_{bc}g^{I_sJ_s}\Psi^c_{J_s}\right)
+\sum_{|\beta|\leq s-1}(Q_\beta)^{aJ_s}_c\partial^\beta\Psi^c_{J_s},
\] with smooth coefficients $Q_\beta$; the fiber indices $b,c$ range over $\{1,\dots,r\}$. The metric factors arise when both inner products are written in components before integrating by parts. Expanding by the Leibniz rule and using $h^{ab}h_{bc}=\delta^a_c$, the term of order $s$ is $(-1)^sg^{I_sJ_s}\partial_{I_s}\Psi_{J_s}^a$. All the remaining terms have order at most $s-1$. Moving $\sqrt{\det(\mathbf{g})}$ inside the leading derivative changes only terms of this lower order. Thus the local expression used below has the indicated leading term and smooth lower-order coefficients $D_\beta$; metric compatibility is precisely the hypothesis that allows us to identify the formal adjoint with the formula in Corollary~\ref{cor: expresion coordenada adjunto formal s derivada covariante}.

\medskip

\noindent\emph{($\Rightarrow$)} First suppose that $\mathbf{F}\restriction_{U}\in W^{m,p}(U,\mathbf{E}\restriction_{U})$. Then, for each $s\le m$, the $s$th weak covariant derivative of $\mathbf{F}\restriction_{U}$ exists; denote it by $\mathbf{G}\in L^{p}\big(U,T^{(0,s)}(TU)\otimes \mathbf{E}\restriction_{U}\big)$. By definition, for every $\boldsymbol{\Psi}\in \Gamma_{c}\big(T^{(0,s)}(TU)\otimes \mathbf{E}\restriction_{U}\big)$ we have \[
\int_{U}\langle \mathbf{F},(\nabla^{s})^{*}\boldsymbol{\Psi}\rangle_{\mathbf{g},\mathbf{h}_{\mathbf{E}}}\,d\lambda_{\mathbf{g}}
=\int_{U}\langle \mathbf{G},\boldsymbol{\Psi}\rangle_{\mathbf{g},\mathbf{h}_{\mathbf{E}}}\,d\lambda_{\mathbf{g}}.
\] In local coordinates, writing $\mathbf{F}= F^{a} \mathbf{e}_a$ and $\boldsymbol{\Psi} = \Psi^{c}_{J_s}\, \mathbf{d}x^{J_s}\otimes \mathbf{e}_c$, this identity reads \[
\int_{U}
h_{bc}\,
F^{b}\,
\big((\nabla^{s})^{*}\boldsymbol{\Psi}\big)^{c}
\,d\lambda_{\mathbf{g}}
=\int_{U}
g^{I_s J_s}\,h_{bc}\,
G^{b}_{I_s}\,
\Psi^{c}_{J_s}
\,d\lambda_{\mathbf{g}}.
\]

Now fix a target block $T_s=(t_{1},\dots,t_{s})$ and a fiber index $a$. Given $\psi\in C_{c}^{\infty}(\phi(U))$, define the test tensor field $\boldsymbol{\Psi}$ by its components: \begin{equation}\label{eq: campo tensorial de prueba conveniente}
\Psi^{c}_{J_s}
=\frac{1}{\sqrt{\det(\mathbf{g})}}\,h^{ca}\,
g_{J_s T_s}\,
(\psi\circ \phi).
\end{equation} The coordinate open set $U$, regarded as a manifold, has no boundary, and the test sections have compact support in it. Thus Theorem~\ref{teo:existencia-adjunto-hermitiano-formal}, in the local form described at the beginning of this proof, applies on $U$, and the formal adjoint acts on $\boldsymbol{\Psi}$ as \[
\big((\nabla^{s})^{*}\boldsymbol{\Psi}\big)^{c} =
(-1)^{s}\left[\frac{1}{\sqrt{\det(\mathbf{g})}}
\partial_{I_s}\left(
\sqrt{\det(\mathbf{g})}g^{I_{s}J_{s}}\Psi^{c}_{J_s}\right)+\displaystyle\sum_{|\beta|\le s-1}
\big(D_{\beta}\big)^{c\, J_s}_{d}
\partial^{\beta}
\Psi^{d}_{J_s}
\right].
 \] Substituting \eqref{eq: campo tensorial de prueba conveniente} into this expression and using $\displaystyle\frac{\sqrt{\det(\mathbf{g})}}{\sqrt{\det(\mathbf{g})}}=1$, $g^{I_{s}J_{s}}g_{J_{s}T_{s}}=\delta^{I_{s}}_{T_{s}}$, we obtain

\[\big((\nabla^{s})^{*}\boldsymbol{\Psi}\big)^{c} =
(-1)^{s}\left[
\frac{1}{\sqrt{\det(\mathbf{g})}}\partial_{T_s}
\left(h^{ca}
(\psi\circ \phi)\right)+\displaystyle\sum_{|\beta|\le s-1}
\big(D_{\beta}\big)^{c\, J_s}_{d}
\partial^{\beta}\left(\frac{h^{da}g_{J_s T_s}}{\sqrt{\det(\mathbf{g})}}
(\psi\circ \phi)\right)\right].\]

Let $I_{1}:=\displaystyle\int_{U}\langle \mathbf{F},(\nabla^{s})^{*}\boldsymbol{\Psi}\rangle_{\mathbf{g},\mathbf{h}_{\mathbf{E}}}\,d\lambda_{\mathbf{g}}$ and $I_{2}:=\displaystyle\int_{U}\langle \mathbf{G},\boldsymbol{\Psi}\rangle_{\mathbf{g},\mathbf{h}_{\mathbf{E}}}\,d\lambda_{\mathbf{g}}$. Substituting \eqref{eq: campo tensorial de prueba conveniente} into $I_{2}$ and noting that $g^{I_s J_s}g_{J_s T_s} = \delta^{I_s}_{T_s}$ gives \begin{equation}\label{eq: I2 ida}
I_{2}=\int_{U}\frac{1}{\sqrt{\det(\mathbf{g})}}\delta^{I_s}_{T_s}\delta_{b}^{a}G^{b}_{I_s}(\psi\circ \phi)\,d\lambda_{\mathbf{g}}
=\int_{\phi(U)}(G^{a}_{T_s}\circ \phi^{-1})\psi \,d\lambda_{n}.
\end{equation}

For $I_{1}$, we use the Leibniz formula in $\partial_{T_s}$:

\[\partial_{T_s}\left(h^{ca}(\psi\circ \phi)\right)
= h^{ca}\partial_{T_s}(\psi\circ \phi) +\displaystyle\sum_{\alpha<\beta_{T_s}}\binom{\beta_{T_s}}{\alpha}\partial^{\alpha}(\psi\circ \phi)\partial^{\beta_{T_s}-\alpha}\left(h^{ca}\right).\]

Substituting back into the integral $I_{1}$ gives \[
I_{1} = (-1)^{s}\int_{U}\frac{1}{\sqrt{\det(\mathbf{g})}}F^{a}\partial_{T_s}(\psi\circ \phi) \,d\lambda_{\mathbf{g}}
\] \[
\quad + (-1)^{s}\displaystyle\sum_{\alpha<\beta_{T_s}}\int_{U}\frac{1}{\sqrt{\det(\mathbf{g})}}h_{bc}F^{b}\binom{\beta_{T_s}}{\alpha}\partial^{\alpha}(\psi\circ \phi)\partial^{\beta_{T_s}-\alpha}\left(h^{ca}\right)\,d\lambda_{\mathbf{g}}
\] \[
\quad + \displaystyle\sum_{|\beta|\le s-1}(-1)^{s}\int_{U}h_{bc}F^{b}\big(D_{\beta}\big)^{c\, J_s}_{d}
\partial^{\beta}\left(\frac{h^{da}g_{J_s T_s}}{\sqrt{\det(\mathbf{g})}}
(\psi\circ \phi)\right)\,d\lambda_{\mathbf{g}}.
\]

We now proceed by induction on $s$. For $s=1$, we have $T_1=(t_1)$, $J_{1}=(j_{1})$, so the expression reduces to \[
I_{1}= -\int_{\phi(U)}(F^{a}\circ \phi^{-1})\partial_{t_{1}}\psi \,d\lambda_{n}
-\int_{\phi(U)}((h_{bc}F^{b})\circ \phi^{-1})\partial_{t_{1}}\left(h^{ca}\circ\phi^{-1}\right)\psi\,d\lambda_{n}
\] \[
\quad -\int_{\phi(U)}\left(\left(h_{bc}F^{b}(D_{0})^{c\, j_1}_{d}h^{da}g_{j_1 t_1}\right)\circ\phi^{-1}\right)\psi\,d\lambda_{n}.\] Define \[H_{t_{1}}^{a}:=h_{bc}F^{b}\partial_{t_{1}}\left(h^{ca}\right)+h_{bc}F^{b}(D_{0})^{c\, j_1}_{d}h^{da}g_{j_1 t_1}.\] Since $\mathbf{F}\restriction_U \in W^{m,p}(U,\mathbf{E}\restriction_{U})$ implies $F^{b}\circ\phi^{-1}\in L^p(\phi(U))$, and the metric coefficients are smooth on the compact set $\overline{\phi(U)}$, it follows that $H_{t_{1}}^{a}\circ \phi^{-1}\in L^{p}(\phi(U))$. Since $I_1 = I_2$, we obtain \[
-\int_{\phi(U)}(F^{a}\circ \phi^{-1})\partial_{t_{1}}\psi \,d\lambda_{n}
=\int_{\phi(U)}\big((G_{t_{1}}^{a}\circ \phi^{-1})+(H^{a}_{t_{1}}\circ \phi^{-1})\big)\psi \,d\lambda_{n},
\] which shows that $F^{a}\circ\phi^{-1}\in W^{1,p}(\phi(U))$.

Now suppose that the assertion holds through order $s-1$ for $s\leq m$. The induction hypothesis ensures that the components of $\mathbf{F}$ on $U$ have weak derivatives of order $\leq s-1$, so we may integrate by parts in the weak sense in the lower-order terms.

First, consider the terms arising from the Leibniz expansion in $I_1$: \[(-1)^{s}\displaystyle\sum_{\alpha<\beta_{T_s}}\int_{U}\frac{h_{bc}F^{b}}{\sqrt{\det(\mathbf{g})}}\binom{\beta_{T_s}}{\alpha}\partial^{\alpha}(\psi\circ \phi)\partial^{\beta_{T_s}-\alpha}\left(h^{ca}\right)\,d\lambda_{\mathbf{g}}\] \[=(-1)^{s}\displaystyle\sum_{\alpha<\beta_{T_s}}\int_{\phi(U)}\left[\left(h_{bc}F^{b}\binom{\beta_{T_s}}{\alpha}\partial^{\beta_{T_s}-\alpha}\left(h^{ca}\right)\right)\circ \phi^{-1}\right]\partial^{\alpha}\psi \,d\lambda_{n}.\]

Since $|\alpha| < |\beta_{I_s}| = s$, we have $|\alpha| \le s-1$. Integrating by parts (in the weak sense) $|\alpha|$ times transfers the derivatives from $\psi$ to the coefficient (which involves $F^b$), so the preceding term equals \[
\displaystyle\sum_{\alpha<\beta_{T_s}}(-1)^{s+|\alpha|}\int_{\phi(U)}\partial^{\alpha}\left[\left(h_{bc}F^{b}\binom{\beta_{T_s}}{\alpha}\partial^{\beta_{T_s}-\alpha}\left(h^{ca}\right)\right)\circ \phi^{-1}\right]
\psi\, d\lambda_{n}.
\]

Second, consider the terms in $I_1$ involving the connection coefficients ($D_{\beta}$). Using the same change of variables and recalling that $|\beta|\le s-1$, we obtain \[
 \displaystyle\sum_{|\beta|\le s-1}(-1)^{s}\int_{U}h_{bc}F^{b}\big(D_{\beta}\big)^{c\, J_s}_{d}
\partial^{\beta}\left(\frac{h^{da}g_{J_s T_s}}{\sqrt{\det(\mathbf{g})}}
(\psi\circ \phi)\right)\,d\lambda_{\mathbf{g}}
\] \[
=(-1)^{s}\displaystyle\sum_{|\beta|\le s-1}\int_{\phi(U)}\left[\left(h_{bc}F^{b}\big(D_{\beta}\big)^{c\, J_s}_{d}\sqrt{\det(\mathbf{g})}\right)\circ \phi^{-1}\right]
\partial^{\beta}\left(\left(\frac{h^{da}g_{J_s T_s}}{\sqrt{\det(\mathbf{g})}}
\circ \phi^{-1}\right)\psi\right)\,d\lambda_{n}.
\] Integrating by parts $|\beta|$ times gives

\[= (-1)^{s}\displaystyle\sum_{|\beta|\le s-1}(-1)^{|\beta|}\int_{\phi(U)}\partial^{\beta}\left[\left(h_{bc}F^{b}\big(D_{\beta}\big)^{c\, J_s}_{d}\sqrt{\det(\mathbf{g})}\right)\circ \phi^{-1}\right]\left(\frac{h^{da}g_{J_s T_s}}{\sqrt{\det(\mathbf{g})}}
\circ \phi^{-1}\right)\psi\,d\lambda_{n}.\]

Now define $\widetilde{G}^{a}_{T_s}$ on $\phi(U)$ by grouping $G^{a}_{T_s}\circ \phi^{-1}$ (which comes from the integral $I_2$) and moving to the right-hand side the terms computed above from $I_1$ (introducing an extra factor $-1$, which changes to $(-1)^{s+1}$): \[
\widetilde{G}^{a}_{T_s} := G^{a}_{T_s}\circ \phi^{-1}
+\displaystyle\sum_{\alpha<\beta_{T_s}}(-1)^{s+1+|\alpha|}\partial^{\alpha}\left[\left(h_{bc}F^{b}\binom{\beta_{T_s}}{\alpha}\partial^{\beta_{T_s}-\alpha}\left(h^{ca}\right)\right)\circ \phi^{-1}\right]
\] \[
+ \displaystyle\sum_{|\beta|\le s-1}(-1)^{s+1+|\beta|}\partial^{\beta}\left[\left(h_{bc}F^{b}\big(D_{\beta}\big)^{c\, J_s}_{d}\sqrt{\det(\mathbf{g})}\right)\circ \phi^{-1}\right]\left(\frac{h^{da}g_{J_s T_s}}{\sqrt{\det(\mathbf{g})}}
\circ \phi^{-1}\right).
\]

The induction hypothesis states that $F^{b}\circ\phi^{-1}\in W^{s-1,p}(\phi(U))$ for every $b\in \{1,\dots,r\}$. Since $\widetilde{G}_{T_s}^{a}$ consists of smooth bounded functions (defined on $\overline{\phi(U)}$) and weak derivatives of $F^{b}\circ\phi^{-1}$ of order at most $s-1$, we conclude that $\widetilde{G}_{T_s}^{a}\in L^{p}(\phi(U))$. Using the equality $I_{1}=I_{2}$ and the definition of $\widetilde G$ gives \[
(-1)^{s}\int_{\phi(U)}(F^{a}\circ \phi^{-1})\partial_{T_s}\psi \,d\lambda_{n}
= \int_{\phi(U)}\widetilde{G}_{T_s}^{a}\psi \,d\lambda_{n}.
\] This shows that $F^{a}\circ \phi^{-1}$ has a weak derivative $\partial_{T_s}(F^a \circ \phi^{-1}) = \widetilde{G}_{T_s}^{a} \in L^{p}(\phi(U))$. Therefore, $F^{a}\circ \phi^{-1}\in W^{s,p}(\phi(U))$.

Since $s\le m$ was arbitrary, we conclude that $F^{a}\circ\phi^{-1}\in W^{m,p}(\phi(U))$ for every $a \in \{1,\dots,r\}$.

\medskip

\noindent \textit{($\Leftarrow$)} Now suppose that $F^{a}\circ\phi^{-1}\in W^{m,p}(\phi(U))$ for every $a\in\{1,\dots,r\}$. Fix $s\le m$ and take $\boldsymbol{\Psi}\in\Gamma_{c}\big(T^{(0,s)}(TU)\otimes \mathbf{E}\restriction_{U}\big)$. We have \[
\int_{U}\langle \mathbf{F},(\nabla^{s})^{*}\boldsymbol{\Psi}\rangle_{\mathbf{g},\mathbf{h}_{\mathbf{E}}}\,d\lambda_{\mathbf{g}}
=
\int_{U}
h_{ab}\,
F^{a}
\big((\nabla^{s})^{*}\boldsymbol{\Psi}\big)^{b}\,
d\lambda_{\mathbf{g}}.
\] By the local expression for the formal adjoint justified at the beginning of this proof, the adjoint has the coordinate form \[
\big((\nabla^{s})^{*}\boldsymbol{\Psi}\big)^{b} = (-1)^{s}\left[\frac{1}{\sqrt{\det(\mathbf{g})}}
\partial_{I_s}\left(
\sqrt{\det(\mathbf{g})}g^{I_{s}J_{s}}\Psi^{b}_{J_s}\right)
+\displaystyle\sum_{|\beta|\le s-1}
\big(D_{\beta}\big)^{b\, J_s}_{c}\,
\partial^{\beta}
\Psi^{c}_{J_s}
\right].
\] The integral can then be written as a sum of two terms, $I_{1}$ and $I_{2}$. Expanding in coordinates gives \[\int_{U}\langle \mathbf{F},(\nabla^{s})^{*}\boldsymbol{\Psi}\rangle_{\mathbf{g},\mathbf{h}_{\mathbf{E}}}\,d\lambda_{\mathbf{g}} = \underbrace{(-1)^{s}\int_{U}h_{ab}F^{a}\frac{1}{\sqrt{\det(\mathbf{g})}}
\partial_{I_s}\left(
\sqrt{\det(\mathbf{g})}g^{I_{s}J_{s}}\Psi^{b}_{J_s}\right)\,d\lambda_{\mathbf{g}}}_{I_{1}}\] \[
 + \underbrace{(-1)^{s}\displaystyle\sum_{|\beta|\leq s-1}\int_{U}h_{ab}F^{a}\big(D_{\beta}\big)^{b\, J_s}_{c}\,
\partial^{\beta}
\Psi^{c}_{J_s}\,d\lambda_{\mathbf{g}}}_{I_{2}}.\] Changing the domain of integration to $\phi(U)$ and writing the compositions explicitly gives \[
I_1 = (-1)^{s}\int_{\phi(U)}\left[\left(h_{ab}F^{a}\right)\circ\phi^{-1}\right]
\partial_{I_s}\left(\left(
\sqrt{\det(\mathbf{g})}g^{I_{s}J_{s}}\Psi^{b}_{J_s}\right)\circ\phi^{-1}\right)\,d\lambda_{n},
\] \[
I_2 = (-1)^{s}\displaystyle\sum_{|\beta|\leq s-1}\int_{\phi(U)}\left[\left(h_{ab}F^{a}(D_{\beta})^{b\, J_s}_{c}\sqrt{\det(\mathbf{g})}\right)\circ \phi^{-1}\right]\partial^{\beta}(\Psi^{c}_{J_s}\circ \phi^{-1})\,d\lambda_{n}.
\]

We now integrate by parts in the weak sense in $I_1$ (which is valid because $F^{a}\circ \phi^{-1}\in W^{m,p}(\phi(U))$ and $s \le m$): \[
I_{1} = \int_{\phi(U)}\partial_{I_s}\left[\left(h_{ab}\circ\phi^{-1}\right)(F^{a}\circ \phi^{-1})\right]\left(\sqrt{\det(\mathbf{g})}g^{I_{s}J_{s}}\Psi^{b}_{J_s}\circ\phi^{-1}\right)\,d\lambda_{n}.
\] Applying the Leibniz rule to the bracketed term gives \[
I_{1} = \int_{\phi(U)}\left[\left(h_{ab}g^{I_s J_s}\sqrt{\det(\mathbf{g})}\circ\phi^{-1}\right)\partial_{I_s}(F^{a}\circ\phi^{-1})\right](\Psi^{b}_{J_s}\circ\phi^{-1})\,d\lambda_{n}
\] \[
\quad + \int_{\phi(U)}\left[\displaystyle\sum_{\substack{I_{s}\\\alpha<\beta_{I_s}}}\binom{\beta_{I_s}}{\alpha}\partial^{\beta_{I_s}-\alpha}\left(h_{ab}\circ\phi^{-1}\right)\partial^{\alpha}(F^{a}\circ \phi^{-1})\left(\left(g^{I_s J_s}\sqrt{\det(\mathbf{g})}\Psi^{b}_{J_s}\right)\circ\phi^{-1}\right)\right]\,d\lambda_{n}.
\] Define the coefficient functions $(B_{\gamma})_{ abI_{s}}$ on $\phi(U)$ for each $|\gamma| \le s-1$ by \[
(B_{\gamma})_{ ab}^{J_s} := \displaystyle\sum_{\substack{I_{s}\\\gamma < \beta_{I_s}}}\binom{\beta_{I_s}}{\gamma}\partial^{\beta_{I_s}-\gamma}\left(h_{ab}\circ\phi^{-1}\right)\left(\left(g^{I_s J_s}\sqrt{\det(\mathbf{g})}\right)\circ\phi^{-1}\right).
\] Substituting this definition into $I_1$ gives \[
I_{1} = \int_{\phi(U)}\left[\left(h_{ab}g^{I_s J_s}\sqrt{\det(\mathbf{g})}\circ\phi^{-1}\right)\partial_{I_s}(F^{a}\circ\phi^{-1})+ \displaystyle\sum_{|\gamma|\leq s-1}(B_{\gamma})_{ ab}^{J_s}\partial^{\gamma}(F^{a}\circ \phi^{-1})\right](\Psi^{b}_{J_s}\circ\phi^{-1})\,d\lambda_{n}.
\]

For $I_{2}$, integrating by parts $|\beta|$ times gives \[
I_{2}=\displaystyle\sum_{|\beta|\le s-1}(-1)^{s+|\beta|}\int_{\phi(U)}\partial^{\beta}\left[\left(\left(h_{ab}(D_{\beta})^{b\, J_s}_{c}\sqrt{\det(\mathbf{g})}\right)\circ\phi^{-1} \right)(F^a \circ \phi^{-1})\right](\Psi^{c}_{J_s}\circ \phi^{-1})\,d\lambda_{n}.
\] Applying the Leibniz rule once more to expand the derivative $\partial^{\beta}$ of the product gives \[
I_{2} = \int_{\phi(U)}\displaystyle\sum_{|\beta|\le s-1}(-1)^{s+|\beta|}\displaystyle\sum_{\gamma \le \beta}\binom{\beta}{\gamma}\partial^{\beta-\gamma}\left(\left(h_{ab}(D_{\beta})^{b\, J_s}_{c}\sqrt{\det(\mathbf{g})}\right)\circ\phi^{-1}\right)\partial^{\gamma}(F^{a}\circ \phi^{-1})(\Psi^{c}_{J_s}\circ \phi^{-1})\,d\lambda_{n}.
\] Define the coefficient functions $(C_{\gamma})^{J_s}_{ac}$ on $\phi(U)$ for each $|\gamma| \le s-1$ by \[
(C_{\gamma})^{J_s}_{ac} := \displaystyle\sum_{\substack{|\beta|\le s-1 \\ \gamma \le \beta}} (-1)^{s+|\beta|} \binom{\beta}{\gamma} \partial^{\beta-\gamma}\left( \left(h_{ab}(D_{\beta})^{b\, J_s}_{c}\sqrt{\det(\mathbf{g})}\right)\circ\phi^{-1} \right).
\] Substituting this definition into $I_2$ gives \[
I_{2} = \int_{\phi(U)}\displaystyle\sum_{|\gamma|\leq s-1}(C_{\gamma})^{J_s}_{ac}\partial^{\gamma}(F^{a}\circ \phi^{-1})(\Psi^{c}_{J_s}\circ \phi^{-1})\,d\lambda_{n}.
\]

We now combine both results. To unify the notation, in $I_1$ we rename the fiber index $b\to c$ (and use $h_{ac}$). In this way, \[\int_{U}\langle \mathbf{F},(\nabla^{s})^{*}\boldsymbol{\Psi}\rangle\,d\lambda _{\mathbf{g}} =\] \[\int_{\phi(U)}\Bigg[\left(h_{ac}g^{I_s J_s}\sqrt{\det(\mathbf{g})}\circ\phi^{-1}\right)\partial_{I_s}(F^{a}\circ\phi^{-1}) + \displaystyle\sum_{|\gamma|\leq s-1}\left((B_{\gamma})_{ ac}^{J_s} + (C_{\gamma})^{J_s}_{ac}\right)\partial^{\gamma}(F^{a}\circ \phi^{-1})\Bigg](\Psi^{c}_{J_s}\circ\phi^{-1})\,d\lambda_{n}.\]

Now define the correction coefficients $(\widehat{D}_{\gamma})$ on $U$ by taking the \textit{pullback} of $B$ and $C$ and inverting the metric weights: \[
\big(\widehat{D}_{\gamma}\big)^{a}_{b\,I_s}
:=\frac{
g_{I_s \tilde{J}_s}h^{a\tilde{c}}
}{\sqrt{\det(\mathbf{g})}}\Big[\big((B_{\gamma})^{\tilde{J}_s}_{b\tilde{c}}\circ\phi\big)+
\big((C_{\gamma})^{\tilde{J}_s}_{b\tilde{c}}
\circ\phi\big)
\Big].
\] Then define $\mathbf{G}\in L^{p}(U,T^{(0,s)}(TU)\otimes \mathbf{E}\restriction_{U})$ in coordinates by \begin{equation}\label{eq: expresion explicita s derivada covariante debil}
G^{a}_{I_s}
=
\partial_{I_s} F^{a}+
\displaystyle\sum_{|\gamma|\le s-1}
\big(\widehat{D}_{\gamma}\big)^{a}_{b\,I_s}
\partial^{\gamma}F^{b}.\end{equation} Let us show that $\mathbf{G}$ is the desired weak covariant derivative of order $s$: \[
\int_{U}\langle \mathbf{G},\boldsymbol{\Psi}\rangle_{\mathbf{g},\mathbf{h}_{\mathbf{E}}}\,d\lambda_{\mathbf{g}} = \int_{U}h_{ac}g^{I_s J_s}G_{I_s}^{a}\Psi_{J_s}^{c}\,d\lambda_{\mathbf{g}}
\] \[
= \int_{U}h_{ac}g^{I_s J_s}\left(\partial_{I_s}F^{a}+\displaystyle\sum_{|\gamma|\leq s-1}(\widehat{D}_{\gamma})^{a}_{b\, I_s}\partial^{\gamma}F^{b}\right)\Psi^{c}_{J_s}\,d\lambda_{\mathbf{g}}.
\] The first term agrees with the principal part. For the lower-order terms, we substitute the definition of $\widehat{D}_{\gamma}$: \[
\int_{U}\displaystyle\sum_{|\gamma|\leq s-1} h_{ac}g^{I_s J_s}\left( \frac{g_{I_s \tilde{J}_s}h^{a\tilde{c}}}{\sqrt{\det(\mathbf{g})}}\Big[\big((B_{\gamma})^{\tilde{J}_s}_{b\tilde{c}}\circ\phi\big)+ \big((C_{\gamma})^{\tilde{J}_s}_{b\tilde{c}}\circ\phi\big)\Big] \right) \partial^{\gamma}F^{b}\Psi^{c}_{J_s} \,d\lambda_{\mathbf{g}}
\] \[
= \int_{U}\displaystyle\sum_{|\gamma|\leq s-1} \frac{1}{\sqrt{\det(\mathbf{g})}} (h_{ac}h^{a\tilde{c}}) (g^{I_s J_s}g_{I_s \tilde{J}_s}) \Big[\big((B_{\gamma})^{\tilde{J}_s}_{b\tilde{c}}\circ\phi\big)+ \big((C_{\gamma})^{\tilde{J}_s}_{b\tilde{c}}\circ\phi\big)\Big] \partial^{\gamma}F^{b}\Psi^{c}_{J_s} \,d\lambda_{\mathbf{g}}.
\] Using the identities $h_{ac}h^{a\tilde{c}} = \delta^{\tilde{c}}_{c}$ and $g^{I_s J_s}g_{I_s \tilde{J}_s} = \delta^{J_s}_{\tilde{J}_s}$ gives \[
= \int_{U}\displaystyle\sum_{|\gamma|\leq s-1} \frac{1}{\sqrt{\det(\mathbf{g})}} \Big[\big((B_{\gamma})^{J_s}_{bc}\circ\phi\big)+ \big((C_{\gamma})^{J_s}_{bc}\circ\phi\big)\Big] \partial^{\gamma}F^{b}\Psi^{c}_{J_s}\,d\lambda_{\mathbf{g}}.
\]

\[
= \int_{\phi(U)}\displaystyle\sum_{|\gamma|\leq s-1} \left[ (B_{\gamma})^{J_s}_{bc} + (C_{\gamma})^{J_s}_{bc} \right] \partial^{\gamma}(F^{b}\circ \phi^{-1})(\Psi^{c}_{J_s}\circ\phi^{-1})\,d\lambda_{n}.\]

This expression agrees exactly with the expansion of $\displaystyle\int_{U}\langle \mathbf{F},(\nabla^{s})^{*}\boldsymbol{\Psi}\rangle\,d\lambda _{\mathbf{g}}$. Therefore, $\mathbf{G}$ is the $s$th weak covariant derivative of $\mathbf{F}$ on $U$. Moreover, $\mathbf{G}\in L^{p}(U,T^{(0,s)}(TU)\otimes \mathbf{E}\restriction_{U})$, since each coordinate component is a linear combination of weak derivatives of the components of $\mathbf{F}$ (through order $s$) multiplied by smooth bounded functions (bounded on $U$ because $U$ is a regular coordinate ball). Since $s\le m$ was arbitrary, we conclude that $\mathbf{F}\restriction_{U}\in W^{m,p}\big(U,\mathbf{E}\restriction_{U}\big)$.

Finally, we prove both norm inequalities quantitatively. Set $n:=\dim M$, $\Omega:=\phi(U)$, and $f^a:=F^a\circ\phi^{-1}$ for $a\in\{1,\dots,r\}$. Recall that \[
\|\mathbf F\|_{W^{m,p}(U,\mathbf E\restriction_U)}
=\left(\sum_{s=0}^{m}\int_U
|\nabla_w^s\mathbf F|_{\mathbf g,\mathbf h_{\mathbf E}}^p
\,d\lambda_{\mathbf g}\right)^{1/p}.
\] At order zero we understand $\nabla_w^0\mathbf F=\mathbf F$.

Define the bundle metric $\overline{\mathbf h}_{\mathbf E}$ by making the frame $\mathbf e_1,\dots,\mathbf e_r$ orthonormal. Explicitly, if $v=\displaystyle\sum_{a=1}^{r}v^a\mathbf e_a(x)$ and $w=\displaystyle\sum_{a=1}^{r}w^a\mathbf e_a(x)$, set \[
\overline{\mathbf h}_{\mathbf E,x}(v,w)
=\sum_{a=1}^{r}v^aw^a.
\] Also let $\delta:=\phi^*\overline{\mathbf g}$, where $\overline{\mathbf g}$ is the Euclidean metric on $\Omega$. The chart and frame are defined on a neighborhood of $\overline U$, so these auxiliary metrics are as well. For each $q\in\{0,\dots,m\}$ define \[
\|\mathbf F\|_{\mathrm{aux},q,p}
:=\left(\sum_{s=0}^{q}\int_U
|\nabla_w^s\mathbf F|_{\delta,\overline{\mathbf h}_{\mathbf E}}^p
\,d\lambda_\delta\right)^{1/p}.
\] The covariant derivatives in this expression remain the original ones: only the fiber norm of each tensor and the measure have changed.

For each $s\in\{0,\dots,m\}$ we apply Lemma~\ref{lema equivalencia metricas fibradas} to the two induced metrics on $T^{(0,s)}(TM)\otimes\mathbf E$ over the compact set $\overline U$. There are numbers $a_s,b_s>0$ such that, for every $x\in\overline U$ and every $\mathbf T\in T^{(0,s)}_x(TM)\otimes\mathbf E_x$, \[
a_s|\mathbf T|_{\delta,\overline{\mathbf h}_{\mathbf E}}
\leq|\mathbf T|_{\mathbf g,\mathbf h_{\mathbf E}}
\leq b_s|\mathbf T|_{\delta,\overline{\mathbf h}_{\mathbf E}}.
\] Taking $a_*:=\min_{0\leq s\leq m}a_s>0$ and $b_*:=\displaystyle\max_{0\leq s\leq m}b_s<\infty$ gives constants valid simultaneously for all orders under consideration. On the other hand, in the coordinates of $\phi$, \[
d\lambda_{\mathbf g}=\rho\,d\lambda_\delta,
\qquad \rho:=\sqrt{\det(\mathbf g)}.
\] The function $\rho$ is smooth and strictly positive on a neighborhood of $\overline U$. By compactness there are $\rho_-,\rho_+>0$ such that $\rho_-\leq\rho(x)\leq\rho_+$ for every $x\in\overline U$. Thus, for each $q\in\{0,\dots,m\}$, \[
\begin{aligned}
a_*^p\rho_-\sum_{s=0}^{q}\int_U
|\nabla_w^s\mathbf F|_{\delta,\overline{\mathbf h}_{\mathbf E}}^p
\,d\lambda_\delta
&\leq
\sum_{s=0}^{q}\int_U
|\nabla_w^s\mathbf F|_{\mathbf g,\mathbf h_{\mathbf E}}^p
\,d\lambda_{\mathbf g}\\
&\leq
b_*^p\rho_+\sum_{s=0}^{q}\int_U
|\nabla_w^s\mathbf F|_{\delta,\overline{\mathbf h}_{\mathbf E}}^p
\,d\lambda_\delta.
\end{aligned}
\] With $c:=a_*\rho_-^{1/p}$ and $C:=b_*\rho_+^{1/p}$, taking $p$th roots and specializing to $q=m$ gives \begin{equation}\label{eq: desigualdad prueba equivalencia normas sobolev localmente}
c\|\mathbf F\|_{\mathrm{aux},m,p}
\leq\|\mathbf F\|_{W^{m,p}(U,\mathbf E\restriction_U)}
\leq C\|\mathbf F\|_{\mathrm{aux},m,p}.
\end{equation} The same comparison holds at each order $q\leq m$ with these constants.

We now pass from covariant derivatives to partial derivatives. For $s\in\{1,\dots,m\}$ write $\mathcal I_s:=\{1,\dots,n\}^s$; its elements are ordered blocks $I_s=(i_1,\dots,i_s)$. For $s=0$ set $\mathcal I_0:=\{\varnothing\}$ and $\partial_{\varnothing}f^a:=f^a$. Formula \eqref{eq: expresion explicita s derivada covariante debil} reads \begin{equation}\label{eq:componentes-derivada-covariante-local}
\bigl((\nabla_w^s\mathbf F)^a_{I_s}\bigr)\circ\phi^{-1}
=\partial_{I_s}f^a+
\sum_{b=1}^{r}\sum_{\substack{\gamma\in\mathbb N_0^n\\|\gamma|\leq s-1}}
A^{abI_s}_{s,\gamma}\,\partial^\gamma f^b,
\end{equation} where \[
A^{abI_s}_{s,\gamma}
:=\bigl((\widehat D_\gamma)^a_{b\,I_s}\bigr)\circ\phi^{-1}.
\] Each coefficient extends smoothly to a neighborhood of $\overline\Omega$ and, in particular, belongs to $L^\infty(\Omega)$. If $m\geq1$, define \[
K:=\max_{\substack{1\leq s\leq m,\ 1\leq a,b\leq r\\
I_s\in\mathcal I_s,\ \gamma\in\mathbb N_0^n,\ |\gamma|\leq s-1}}
\|A^{abI_s}_{s,\gamma}\|_{L^\infty(\Omega)}<\infty.
\] The maximum contains finitely many coefficients. If $m=0$, set $K:=0$, since there are no lower-order terms.

Denote by \[
G^a_{s,I_s}:=\bigl((\nabla_w^s\mathbf F)^a_{I_s}\bigr)\circ\phi^{-1},
\qquad
Q_s:=\left(\sum_{a=1}^{r}\sum_{I_s\in\mathcal I_s}
|G^a_{s,I_s}|^2\right)^{1/2}.
\] In particular, $G^a_{0,\varnothing}=f^a$. The coordinate frame tensors, together with the frame $\mathbf e_a$, are orthonormal for the auxiliary metrics. Moreover, $\phi_*\lambda_\delta=\lambda_n$. These two observations give the exact equalities \[
\begin{aligned}
\|\mathbf F\|_{\mathrm{aux},m,p}
&=\left(\sum_{s=0}^{m}\int_\Omega Q_s(x)^p\,d\lambda_n(x)\right)^{1/p}\\
&=\left[
\sum_{s=0}^{m}\int_\Omega
\left(\sum_{a=1}^{r}\sum_{I_s\in\mathcal I_s}
\left|\partial_{I_s}f^a+
\sum_{b=1}^{r}\sum_{\substack{\gamma\in\mathbb N_0^n\\|\gamma|\leq s-1}}
A^{abI_s}_{s,\gamma}\partial^\gamma f^b\right|^2
\right)^{p/2}\,d\lambda_n
\right]^{1/p}.
\end{aligned}
\] In the summand $s=0$, the inner sum of lower-order terms is zero, and the expression inside the bars is simply $f^a$.

For the upper bound, fix $s\in\{1,\dots,m\}$. The inequality $|A+B|^2\leq2|A|^2+2|B|^2$ and the bound $K$ on the coefficients imply, almost everywhere on $\Omega$, \[
\begin{aligned}
Q_s^2
&\leq
2\sum_{a=1}^{r}\sum_{I_s\in\mathcal I_s}|\partial_{I_s}f^a|^2\\
&\quad+
2rn^sK^2
\left(\sum_{b=1}^{r}
\sum_{\substack{\gamma\in\mathbb N_0^n\\|\gamma|\leq s-1}}
|\partial^\gamma f^b|\right)^2\\
&=
2\sum_{a=1}^{r}\sum_{\substack{\alpha\in\mathbb N_0^n\\|\alpha|=s}}
\frac{s!}{\alpha!}|\partial^\alpha f^a|^2
 +2rn^sK^2
\left(\sum_{b=1}^{r}
\sum_{\substack{\gamma\in\mathbb N_0^n\\|\gamma|\leq s-1}}
|\partial^\gamma f^b|\right)^2.
\end{aligned}
\] The factor $rn^s$ counts the pairs $(a,I_s)$. The factor $s!/\alpha!$ counts the ordered blocks producing the same derivative $\partial^\alpha$; distributional partial derivatives commute, so this grouping is valid for weak derivatives as well. Since $s!/\alpha!\leq s!$ and the Euclidean norm of a finite family is at most the sum of the absolute values of its components, it follows that \[
Q_s\leq
\sqrt2\bigl(\sqrt{s!}+\sqrt{rn^s}\,K\bigr)
\sum_{b=1}^{r}\sum_{\substack{\beta\in\mathbb N_0^n\\|\beta|\leq m}}
|\partial^\beta f^b|.
\] At order zero the same bound holds with the factor $1$. Set \[
B_0:=1,\qquad
B_s:=\sqrt2\bigl(\sqrt{s!}+\sqrt{rn^s}\,K\bigr)\ (1\leq s\leq m),
\qquad B:=\max_{0\leq s\leq m}B_s.
\] There are $N_m:=\binom{n+m}{m}$ multi-indices of order at most $m$. Applying the finite-sum inequality \[
\left(\sum_{j=1}^{N}t_j\right)^p
\leq N^{p-1}\sum_{j=1}^{N}t_j^p
\qquad(t_j\geq0,\ N\in\mathbb N),
\] with $N=rN_m$, and then using the definition of the Euclidean Sobolev norm, gives \[
\begin{aligned}
\|\mathbf F\|_{\mathrm{aux},m,p}^p
&\leq (m+1)B^p(rN_m)^{p-1}
\sum_{a=1}^{r}
\sum_{\substack{\beta\in\mathbb N_0^n\\|\beta|\leq m}}
\int_\Omega|\partial^\beta f^a|^p\,d\lambda_n\\
&=(m+1)B^p(rN_m)^{p-1}
\sum_{a=1}^{r}\|f^a\|_{W^{m,p}(\Omega)}^p.
\end{aligned}
\] Thus, with $C_{\mathrm{aux},m}:=(m+1)^{1/p}B(rN_m)^{1-1/p}$, taking roots and using $(\displaystyle\sum_{a=1}^{r}t_a^p)^{1/p}\leq\displaystyle\sum_{a=1}^{r}t_a$ yields \begin{equation}\label{eq:cota-local-superior-componentes}
\|\mathbf F\|_{\mathrm{aux},m,p}
\leq C_{\mathrm{aux},m}
\sum_{a=1}^{r}\|f^a\|_{W^{m,p}(\Omega)}.
\end{equation}

For the lower bound, we proceed by induction on the order. For $q\in\{0,\dots,m\}$, write \[
S_q:=\sum_{a=1}^{r}\|f^a\|_{W^{q,p}(\Omega)},
\qquad U_q:=\|\mathbf F\|_{\mathrm{aux},q,p}.
\] The inductive assertion is that there are $K_q>0$, independent of $\mathbf F$, such that \begin{equation}\label{eq:afirmacion-inductiva-cota-local-inferior}
S_q\leq K_qU_q.
\end{equation} For $q=0$, the pointwise inequality $|f^a|\leq Q_0$ gives $\|f^a\|_{L^p(\Omega)}\leq U_0$ for each $a\in\{1,\dots,r\}$. Summing gives $S_0\leq rU_0$, so we may take $K_0:=r$.

Suppose the result has been proved for $q=s-1$, with $1\leq s\leq m$. For each $\alpha\in\mathbb N_0^n$ with $|\alpha|=s$, choose a block $I(\alpha)\in\mathcal I_s$ containing the index $j$ exactly $\alpha_j$ times for each $j\in\{1,\dots,n\}$. Solving in \eqref{eq:componentes-derivada-covariante-local} gives \[
\partial^\alpha f^a
=G^a_{s,I(\alpha)}
-\sum_{b=1}^{r}\sum_{\substack{\gamma\in\mathbb N_0^n\\|\gamma|\leq s-1}}
A^{abI(\alpha)}_{s,\gamma}\partial^\gamma f^b.
\] Since $|G^a_{s,I(\alpha)}|\leq Q_s$, the triangle inequality in $L^p(\Omega)$ and the estimate for bounded multipliers give \[
\|\partial^\alpha f^a\|_{L^p(\Omega)}
\leq U_s+
K\sum_{b=1}^{r}\sum_{\substack{\gamma\in\mathbb N_0^n\\|\gamma|\leq s-1}}
\|\partial^\gamma f^b\|_{L^p(\Omega)}.
\] For each component $b$, Hölder's inequality for the finite sum of its $N_{s-1}=\binom{n+s-1}{s-1}$ derivatives gives \[
\begin{aligned}
\sum_{\substack{\gamma\in\mathbb N_0^n\\|\gamma|\leq s-1}}
\|\partial^\gamma f^b\|_{L^p(\Omega)}
&\leq N_{s-1}^{1-1/p}
\left(\sum_{\substack{\gamma\in\mathbb N_0^n\\|\gamma|\leq s-1}}
\|\partial^\gamma f^b\|_{L^p(\Omega)}^p\right)^{1/p}\\
&=N_{s-1}^{1-1/p}\|f^b\|_{W^{s-1,p}(\Omega)}.
\end{aligned}
\] For $p=1$ this inequality is an equality. There are $d_s:=\binom{n+s-1}{s}$ multi-indices of order exactly $s$. Thus, summing the preceding bound over the $r$ components and these $d_s$ multi-indices gives \[
\sum_{a=1}^{r}\sum_{\substack{\alpha\in\mathbb N_0^n\\|\alpha|=s}}
\|\partial^\alpha f^a\|_{L^p(\Omega)}
\leq rd_sU_s+rd_sK N_{s-1}^{1-1/p}S_{s-1}.
\] The definition of the Sobolev norm also implies \[
\begin{aligned}
S_s
&\leq S_{s-1}
+\sum_{a=1}^{r}\sum_{\substack{\alpha\in\mathbb N_0^n\\|\alpha|=s}}
\|\partial^\alpha f^a\|_{L^p(\Omega)}\\
&\leq rd_sU_s+
\bigl(1+rd_sK N_{s-1}^{1-1/p}\bigr)S_{s-1}.
\end{aligned}
\] By the induction hypothesis, $S_{s-1}\leq K_{s-1}U_{s-1}\leq K_{s-1}U_s$; the second inequality follows because the auxiliary norm of order $s$ contains all the terms of the norm of order $s-1$. Consequently, \eqref{eq:afirmacion-inductiva-cota-local-inferior} holds for $q=s$ with \[
K_s:=rd_s+
\bigl(1+rd_sK N_{s-1}^{1-1/p}\bigr)K_{s-1}.
\] The induction reaches $s=m$ and proves \begin{equation}\label{eq:cota-local-inferior-componentes}
K_m^{-1}\sum_{a=1}^{r}\|f^a\|_{W^{m,p}(\Omega)}
\leq\|\mathbf F\|_{\mathrm{aux},m,p}.
\end{equation}

Finally, combine \eqref{eq:cota-local-superior-componentes}, \eqref{eq:cota-local-inferior-componentes}, and \eqref{eq: desigualdad prueba equivalencia normas sobolev localmente}. Taking $c_m:=c/K_m$ and $C_m:=C C_{\mathrm{aux},m}$ gives \[
\begin{aligned}
c_m\sum_{a=1}^{r}
\|F^a\circ\phi^{-1}\|_{W^{m,p}(\phi(U))}
&\leq\|\mathbf F\|_{W^{m,p}(U,\mathbf E\restriction_U)},\\
\|\mathbf F\|_{W^{m,p}(U,\mathbf E\restriction_U)}
&\leq C_m\sum_{a=1}^{r}
\|F^a\circ\phi^{-1}\|_{W^{m,p}(\phi(U))}.
\end{aligned}
\] The constants depend on $m,p,n,r$, the fixed chart and frame, the comparison bounds for the metrics and volume on $\overline U$, and the $L^\infty(\Omega)$ norms of the coefficients appearing through order $m$. None depends on $\mathbf F$. The argument includes $m=0$ under the indicated conventions and does not use approximation by smooth sections.

\end{proof}

\begin{example}[A Sobolev section with a corner in its components]
\label{ej:seccion-sobolev-esquina-marcos}
Let $(M,\mathbf g)$ be a Riemannian manifold without boundary of
dimension $n\geq1$, and let $\mathbf E\to M$ be a real bundle of rank
$r\geq1$ with a bundle metric and compatible connection.
Choose a regular coordinate ball
$\phi\colon U\to B_{\mathrm{euc}}(0,R)$ and a smooth frame
$(\mathbf e_a)_{a=1}^r$ defined on a neighborhood of $\overline U$.
Take $\chi\in C_c^\infty(B_{\mathrm{euc}}(0,R))$ equal to $1$ near the
origin, and set
\[
 v(x):=\chi(x)|x^1|,
 \qquad
 \boldsymbol\sigma(p):=
 \begin{cases}
  v(\phi(p))\mathbf e_1(p),&p\in U,\\
  0_p,&p\notin U.
 \end{cases}
\]
This section is continuous and has compact support contained in $U$.
It is not $C^1$ near $\phi^{-1}(0)$, since its first component is
$|x^1|$ where $\chi=1$, and the criterion in
Proposition~\ref{suavidad de secciones en bruto locales} detects
the failure of differentiability.

Its components nevertheless belong to $W^{1,p}(\phi(U))$ for every
$1\leq p<\infty$. By Fubini and
Example~\ref{ej:sobolev-potencias-intervalo},
\[
 D_iv=(D_i\chi)|x^1|+
       \chi\,\delta_{i1}\operatorname{sgn}(x^1).
\]
All these functions are bounded and compactly supported.
Lemma~\ref{lema:meyers-serrin-haz-E} gives local membership in the
Sobolev space of sections. The weak covariant derivative and global
membership can also be checked directly.
Write the connection matrices in the frame as
\[
 \nabla^{\mathbf E}_{\boldsymbol\partial_i}\mathbf e_b
 =\sum_{a=1}^r A_i{}^a{}_b\,\mathbf e_a.
\]
Then
\begin{equation}\label{eq:derivada-seccion-esquina}
 \nabla^{\mathbf E}_w\boldsymbol\sigma
 =\sum_{i=1}^n\sum_{a=1}^r
 \bigl(D_iv\,\delta_{a1}+vA_i{}^a{}_1\bigr)
 \,\mathbf d x^i\otimes\mathbf e_a
\end{equation}
on $U$, and it vanishes outside the support.

To justify the weak identity, extend $v$ by zero to $\mathbb R^n$ and
use $v_\varepsilon=v*\rho_\varepsilon$.
The construction in Lemma~\ref{lem:localizar-ejemplos-sobolev-escalar}
gives $v_\varepsilon\to v$ and $D_iv_\varepsilon\to D_iv$ in $L^p$,
with all supports in a fixed compact subset of $\phi(U)$.
The smooth sections
$\boldsymbol\sigma_\varepsilon=(v_\varepsilon\circ\phi)\mathbf e_1$,
extended by zero, satisfy~\eqref{eq:derivada-seccion-esquina} with
$v_\varepsilon$ in place of $v$ by the Leibniz rule.
The matrices $A_i$, frame coefficients, and metrics are bounded on this
compact set; thus the sections and their derivatives converge in the
corresponding $L^p$ spaces.
Passing to the limit in the identity against every compactly supported
test section gives the stated weak derivative.
Hence $\boldsymbol\sigma\in W^{1,p}(M,\mathbf E)$ for every
$1\leq p<\infty$, although it is not a $C^1$ section.

The term $vA_i{}^a{}_1$ explains why differentiating only the components
does not by itself produce a frame-independent object.
If $\widetilde{\mathbf e}=\mathbf e S$, where $S$ is a smooth invertible
matrix, the component columns and connection matrices change according to
\[
 \widetilde s=S^{-1}s,\qquad
 \widetilde A_i=S^{-1}A_iS+S^{-1}(\partial_iS).
\]
The weak product rule and
$\partial_i(S^{-1})=-S^{-1}(\partial_iS)S^{-1}$ give
\[
 D_i\widetilde s+\widetilde A_i\widetilde s
 =S^{-1}(D_is+A_is).
\]
The terms differentiating the change of frame cancel.
This identity verifies that the weak derivative is the same tensor-valued
section in both frames. On a compact set, $S,S^{-1}$ and their first
derivatives are bounded, so both descriptions preserve membership in
$W^{1,p}$. Sobolev regularity is recognized in the components, while the
connection organizes their differentiation intrinsically.
\end{example}

Before formulating the result for sections, we recall the Euclidean version to be used in local trivializations.

\begin{definition}\label{def:Hmp-euclidiano-meyers-serrin}\index{Euclidean Sobolev space by completion@Euclidean Sobolev space by completion} Let $\Omega\subseteq\mathbb R^n$ be open, let $m\in\mathbb N_0$, and let $1\leq p<\infty$. Define $C^{m,p}(\Omega):=\{u\in C^\infty(\Omega)\mid\|u\|_{W^{m,p}(\Omega)}<\infty\}$ and denote its completion with respect to the norm $\|\cdot\|_{W^{m,p}(\Omega)}$ by $H^{m,p}(\Omega)$. \end{definition}

\begin{theorem}[Meyers--Serrin]\label{teorema meyers serrin}\index{Meyers Serrin@Meyers--Serrin} Let $\Omega\subseteq\mathbb R^n$ be open, $m\in\mathbb N_0$, and $1\leq p<\infty$. Then $H^{m,p}(\Omega)=W^{m,p}(\Omega)$. \end{theorem}

\begin{proof}
The inclusion $I_0\colon C^{m,p}(\Omega)\to W^{m,p}(\Omega)$
is an isometry, since both spaces use the same norm.
Theorem~\ref{teo:sobolev-en-abiertos-euclidianos-espacio-banach-2}
gives completeness of the target; for $m=0$ use completeness
of $L^p(\Omega)$ from
Theorem~\ref{teo:b6-espacios-lp-espacios-de-lebesgue}.
If $(f_k)$ is Cauchy in $C^{m,p}(\Omega)$, then
$(I_0f_k)$ converges in $W^{m,p}(\Omega)$. Define
\[
 I([(f_k)]):=\lim_{k\to\infty}I_0f_k.
\]
Two sequences representing the same element of the completion
satisfy $\|f_k-g_k\|_{W^{m,p}}\to0$.
Their images therefore have the same limit, so $I$ is
well defined. Linearity follows by passing to the limit
in sums and scalar multiples, and continuity of the norm gives
\[
 \|I([(f_k)])\|_{W^{m,p}(\Omega)}
 =\lim_{k\to\infty}\|f_k\|_{W^{m,p}(\Omega)}
 =\|[(f_k)]\|_{H^{m,p}(\Omega)}.
\]
In particular, $I$ is an isometry and injective.

For $u\in W^{m,p}(\Omega)$,
Proposition~\ref{prop:densidad-sobolev-orden-entero} gives
$f_k\in C^\infty(\Omega)\cap W^{m,p}(\Omega)$ with
$\|f_k-u\|_{W^{m,p}(\Omega)}<1/k$. Every $f_k$ belongs
to $C^{m,p}(\Omega)$, and
\[
 \|f_k-f_\ell\|_{W^{m,p}(\Omega)}<\frac{1}{k}+\frac{1}{\ell}.
\]
Its class defines an element of $H^{m,p}(\Omega)$ whose
image under $I$ is $u$. Hence $I$ is surjective and gives
the asserted isometric identification.
\end{proof}

\subsection{The Meyers--Serrin approximation theorem in vector bundles} The bundle version of Meyers--Serrin compares the definition using weak covariant derivatives with the completion of smooth sections. To formulate it, we first introduce the space $H^{m,p}(M,\mathbf{E})$, analogous to the construction in Definition~\ref{def: sobolev en variedades riemannianas con completaciones Hmp(M) Wmp(M) funciones escalares}. \begin{definition}\label{def:definicion-sobolev-haces-norma-variedad-riemanniana-haz-vectorial}\index{Sobolev space@Sobolev space!of sections by completion}\index{vector bundle!Sobolev space of sections} Let $(M,\mathbf{g})$ be a Riemannian manifold without boundary, let $\mathbf{E}\longrightarrow M$ be a smooth vector bundle equipped with a bundle metric $\mathbf{h}_{\mathbf{E}}$ and a connection $\nabla^{\mathbf{E}}$, and let $m\in\mathbb N_0$ and $1\leq p<\infty$. Define \[C^{m,p}(M,\mathbf{E}):=\{\mathbf{u}\in \Gamma(\mathbf{E})\mid \|\mathbf{u}\|_{W^{m,p}(M,\mathbf{E})}<\infty \},\] where $\|\mathbf{u}\|_{W^{m,p}(M,\mathbf{E})}:=\left(\displaystyle\sum_{s=0}^{m}\int_{M}|\nabla^{s}\mathbf{u}|_{\mathbf{g},\mathbf{h}_{\mathbf{E}}}^{p}d\lambda_{\mathbf{g}}\right)^{\displaystyle\frac{1}{p}}$.

Define \[H^{m,p}(M,\mathbf{E}):=\overline{C^{m,p}(M,\mathbf{E})}^{\|\cdot\|_{W^{m,p}(M,\mathbf{E})}},\] which is the completion of $C^{m,p}(M,\mathbf{E})$ with respect to the Sobolev norm $\|\cdot\|_{W^{m,p}(M,\mathbf{E})}$. \end{definition}

\begin{example}[A sequence of tangent fields whose limit is not smooth]
\label{ej:incompletitud-secciones-haz-tangente-esfera}
Consider $\mathbf E=TS^2$ over the unit sphere with its induced
metric and Levi--Civita connection. Fix
$m\in\mathbb N_0$ and $1\leq p<\infty$.
Let $h(p)=z$ and let
\[
 \mathbf X(p):=e_3-h(p)p.
\]
Since $\langle e_3-h(p)p,p\rangle=0$, it is a tangent field.
Its ambient expression is smooth and, by the tangent bundle
trivializations, defines a global smooth section. Moreover,
\[
 |\mathbf X(p)|^2=1-h(p)^2\leq1.
\]
In particular, the field does not vanish on the equator.

We can also compute its first covariant derivative.
If $P_p(v)=v-\langle v,p\rangle p$ is orthogonal projection,
the induced connection is $\nabla_Y\mathbf Z=P_p(D_Y\mathbf Z)$.
Indeed, differentiating the ambient inner product proves
metric compatibility, since the fields paired with it are
tangent and projection does not change that inner product.
Its torsion is zero. In a local parametrization $\Phi$, if
$Y=Y^i\partial_i\Phi$ and $\mathbf Z=Z^j\partial_j\Phi$, the
product rule and $\partial_{ij}\Phi=\partial_{ji}\Phi$ give
\[
 D_Y\mathbf Z-D_{\mathbf Z}Y
 =(Y^i\partial_iZ^j-Z^i\partial_iY^j)\partial_j\Phi
 =[Y,\mathbf Z].
\]
This field is tangent and projection preserves it.
Uniqueness of Levi--Civita identifies it with the chosen connection.
Since
\[
 D_Y\mathbf X=-Y(h)p-hY,
\]
we obtain $\nabla_Y\mathbf X=-hY$. For an orthonormal basis
$Y_1,Y_2$ of the tangent space,
\[
 |\nabla\mathbf X|^2
 =\displaystyle\sum_{i=1}^2|\nabla_{Y_i}\mathbf X|^2=2h^2.
\]
All higher derivatives are smooth and bounded by compactness.
Write
\[
 M_j:=\displaystyle\sup_{S^2}|\nabla^j\mathbf X|,\qquad 0\leq j\leq m.
\]

Take the smooth functions $u_k=a_k(h)$ of
Example~\ref{ej:incompletitud-sobolev-esfera-orden-arbitrario},
with $a_k=\rho_{\varepsilon_k}*H_m$, and define
\[
 \mathbf u_k:=u_k\mathbf X\in\Gamma(TS^2),\qquad
 \mathbf U:=H_m(h)\mathbf X.
\]
By compactness every $\mathbf u_k$ has finite Sobolev norm.
We justify control of all its derivatives. For a smooth
function $f$, the product connection's Leibniz rule gives
\[
 |\nabla^s(f\mathbf X)|
 \leq\displaystyle\sum_{j=0}^s\binom{s}{j}
                   |\nabla^jf|\,|\nabla^{s-j}\mathbf X|.
\]
To prove it, at each of the $s$ differentiations choose
whether to differentiate the scalar or the vector factor.
There are $\binom{s}{j}$ choices with $j$ differentiations
of the first factor. Each term preserves the order of the
slots of each factor and only permutes their slots when
combining them. The product connection commutes with these
permutations, which are isometries for the tensor metric.
The norm of each tensor product is the product of its factors'
norms. The triangle inequality proves the pointwise bound,
including $s=0$. Taking $L^p$ norms and summing the orders,
\[
 \begin{aligned}
 \|f\mathbf X\|_{W^{m,p}(S^2,TS^2)}
 &\leq\displaystyle\sum_{s=0}^m\displaystyle\sum_{j=0}^s
           \binom{s}{j}M_{s-j}\|\nabla^jf\|_{L^p(S^2)}\\
 &\leq C_m\|f\|_{W^{m,p}(S^2)},\qquad
 C_m:=\displaystyle\sum_{s=0}^m\displaystyle\sum_{j=0}^s\binom{s}{j}M_{s-j}<\infty.
 \end{aligned}
\]
The first inequality uses that the norm given by the sum
of powers is at most the sum of the individual orders'
norms; the second uses that each order's norm is at most
the full Sobolev norm. Applying the estimate to $f=u_k-u_\ell$,
we obtain
\[
 \|\mathbf u_k-\mathbf u_\ell\|_{W^{m,p}(S^2,TS^2)}
 \leq C_m\|u_k-u_\ell\|_{W^{m,p}(S^2)}
 \leq C_m(c_k+c_\ell)\longrightarrow0.
\]
The sequence is Cauchy. Completeness of the space defined
by weak covariant derivatives, proved in
Theorem~\ref{encaje isometrico en Lp sobolev haces}, gives its limit.
Since $|\mathbf X|\leq1$,
\[
 \|\mathbf u_k-\mathbf U\|_{L^p(S^2,TS^2)}
 \leq\|u_k-H_m(h)\|_{L^p(S^2)}\longrightarrow0.
\]
The Sobolev norm controls the order-zero norm, so uniqueness
of the $L^p$ limit identifies the Sobolev limit with the
class of $\mathbf U$. In particular,
$\mathbf U\in W^{m,p}(S^2,TS^2)$.

Suppose some smooth section $\mathbf V$ represents this same
class. If $m\geq1$, $\mathbf U$ is continuous, and equality
almost everywhere implies equality at every point: if their
values differed at a point, the continuous norm of the
difference would be positive on an open set of positive measure.
Along the curve $\gamma(t)=(\sqrt{1-t^2},0,t)$, $|t|<1/2$,
the third ambient component satisfies
\[
 \mathbf U_3(\gamma(t))=H_m(t)(1-t^2).
\]
Composition of a smooth section with this curve, followed
by the ambient functional $v\mapsto v_3$, is smooth.
Since $1-t^2>0$, equality $\mathbf V=\mathbf U$ would require
\[
 \frac{\mathbf V_3(\gamma(t))}{1-t^2}=H_m(t)
\]
to be smooth. However, the derivative of order $m$ of the
right-hand side is zero for $t<0$ and one for $t>0$, so it
is not $C^m$ at zero. If $m=0$, continuity forces $\mathbf V=0$
on the negative hemisphere and $\mathbf V=\mathbf X$ on the
positive hemisphere. The limits of $\mathbf V_3(\gamma(t))$
at zero would then be zero from the left and one from the
right, another contradiction. Thus $\mathbf U$ has no smooth
representative, and $C^{m,p}(S^2,TS^2)$ is incomplete in the
Sobolev norm.
\end{example}
 As in the Euclidean case of Theorem~\ref{teorema meyers serrin}, both constructions agree. The equality $H=W$ thus extends from functions to sections of a vector bundle. \begin{theorem}[Meyers--Serrin for vector bundles]\label{meyers-serrin-haz}\index{Meyers Serrin for vector bundles@Meyers--Serrin for vector bundles} Let $(M,\mathbf{g})$ be a Riemannian manifold without boundary, let $\pi_{\mathbf{E}}\colon \mathbf{E} \longrightarrow M$ be a smooth vector bundle equipped with a bundle metric $\mathbf{h}_{\mathbf{E}}$ and a compatible connection $\nabla^{\mathbf{E}}$, and let $m\in\mathbb N_0$ and $1\leq p<\infty$. Then \[
H^{m,p}(M,\mathbf{E})=W^{m,p}(M,\mathbf{E}).
\] \end{theorem}

\begin{proof} By Proposition~\ref{bolacoordenadaregular}, there is a countable locally finite open cover of $M$ by regular coordinate balls $(U_{\alpha},\phi_{\alpha})_{\alpha\in\mathbb{N}}$. Refining it if necessary, we may require that $\overline{U_\alpha}$ lie in the domain of a local frame of $\mathbf{E}$. Let $(\psi_{\alpha})_{\alpha\in\mathbb{N}}$ be a smooth partition of unity subordinate to this cover.

Now consider a section $\mathbf{F}\in W^{m,p}(M,\mathbf{E})$. We can write \[
\mathbf{F}=\displaystyle\sum_{\alpha=1}^{\infty}\psi_{\alpha}\mathbf{F},
\] where the sum is locally finite. For each $\alpha$, the product $\psi_{\alpha}\mathbf{F}\restriction_{U_{\alpha}}$ belongs to $W^{m,p}(U_{\alpha},\mathbf{E}\restriction_{U_{\alpha}})$, since $\psi_{\alpha}$ is smooth and has compact support contained in $U_{\alpha}$.

Fix $\alpha$ and work in the chart $(U_{\alpha},\phi_{\alpha})$, together with one of these local frames $(\mathbf{e}_{1},\dots,\mathbf{e}_{r})$, defined on a neighborhood of $\overline{U_\alpha}$. Applying Lemma~\ref{lema:meyers-serrin-haz-E} to $(\psi_{\alpha}\mathbf{F})\restriction_{U_{\alpha}}$ gives, for every fiber index $a \in \{1,\dots,r\}$, \[
\big(\psi_{\alpha}\mathbf{F}\big)^{a}\circ \phi_{\alpha}^{-1}\in W^{m,p}(\phi_{\alpha}(U_{\alpha})).
\] Moreover, there are constants $C_{m,\alpha}>0$ such that, for every section $\mathbf{H}\in W^{m,p}(U_{\alpha},\mathbf{E}\restriction_{U_{\alpha}})$, \[
\|\mathbf{H}\|_{W^{m,p}(U_{\alpha},\mathbf{E}\restriction_{U_{\alpha}})}\leq
C_{m,\alpha}\displaystyle\sum_{a=1}^{r}
\big\|H^{a}\circ \phi_{\alpha}^{-1}\big\|_{W^{m,p}(\phi_{\alpha}(U_{\alpha}))}.
\]

Write $\Omega_{\alpha}:=\phi_{\alpha}(U_{\alpha})\subset\mathbb{R}^{n}$. The support of each component $\big(\psi_\alpha \mathbf{F}\big)^a\circ\phi_\alpha^{-1}$ is a compact subset of $\Omega_\alpha$. The classical Meyers--Serrin theorem on domains in $\mathbb{R}^{n}$, Theorem~\ref{teorema meyers serrin}, provides smooth approximants in $W^{m,p}(\Omega_\alpha)$. To obtain compactly supported approximants, choose $\chi_\alpha\in C_c^\infty(\Omega_\alpha)$ equal to one on a neighborhood of the preceding support and multiply by $\chi_\alpha$. The Leibniz rule shows directly that multiplication by $\chi_\alpha$ is continuous on $W^{m,p}(\Omega_\alpha)$. Thus, for each $\alpha\in \mathbb{N}$ and each component $a\in\{1,\dots,r\}$, there is a sequence \[
(\widetilde G_{\alpha,k})^{a}\in C_{c}^{\infty}(\Omega_{\alpha}),\qquad k\in \mathbb{N},
\] such that \[
\displaystyle\sum_{a=1}^{r}
\left\|
(\widetilde G_{\alpha,k})^{a}-
\big(\psi_{\alpha}\mathbf{F}\big)^{a}\circ \phi_{\alpha}^{-1}
\right\|_{W^{m,p}(\Omega_{\alpha})}
\xrightarrow[k\to\infty]{}0.
\]

From these smooth functions we define a smooth section $\mathbf{G}_{\alpha,k}\in \Gamma(\mathbf{E}\restriction_{U_{\alpha}
})$ by \[
\mathbf{G}_{\alpha,k}
=
\displaystyle\sum_{a=1}^{r}
\Big[(\widetilde G_{\alpha,k})^{a}\circ \phi_{\alpha}\Big]\; \mathbf{e}_a.
\] The estimate in Lemma~\ref{lema:meyers-serrin-haz-E} then gives \[
\|\mathbf{G}_{\alpha,k}-(\psi_{\alpha}\mathbf{F})\restriction_{U_{\alpha}}\|_{W^{m,p}(U_{\alpha},\mathbf{E}\restriction_{U_{\alpha}})}\xrightarrow[k\to\infty]{}0.
\]

Extending each $\mathbf{G}_{\alpha,k}$ by zero outside $U_{\alpha}$, we may regard it as a global compactly supported section, that is, $\mathbf{G}_{\alpha,k}\in \Gamma_{c}(\mathbf{E})$, and moreover \[\|\mathbf{G}_{\alpha,k}-(\psi_{\alpha}\mathbf{F})\restriction_{U_{\alpha}}\|_{W^{m,p}(U_{\alpha},\mathbf{E}\restriction_{U_{\alpha}})}=\left(\displaystyle\sum_{s=0}^{m}\int_{U_{\alpha}}|\nabla_{w}^{s}(\mathbf{G}_{\alpha,k}-(\psi_{\alpha} \mathbf{F})\restriction_{U_{\alpha}})|^{p}_{\mathbf{g},\mathbf{h}_{\mathbf{E}}}d\lambda_{\mathbf{g}}\right)^{\frac{1}{p}}.\] Since $\supp(\mathbf{G}_{\alpha,k}-\psi_{\alpha}\mathbf{F})\subseteq U_{\alpha}$ and $\nabla^{s}_{w}$ is a local operator (Proposition~\ref{prop: derivada covariante debil es operador local}), the preceding expression equals \[\left(\displaystyle\sum_{s=0}^{m}\int_{M}|\nabla_{w}^{s}(\mathbf{G}_{\alpha,k}-\psi_{\alpha} \mathbf{F})|^{p}_{\mathbf{g},\mathbf{h}_{\mathbf{E}}}d\lambda_{\mathbf{g}}\right)^{\frac{1}{p}}.\] Consequently,

\[
\|\mathbf{G}_{\alpha,k}-\psi_{\alpha}\mathbf{F}\|_{W^{m,p}(M,\mathbf{E})}\xrightarrow[k\to\infty]{}0.
\]

Now fix $\alpha\in \mathbb{N}$ and an integer $q\in \mathbb{N}$. By the preceding convergence, there is $k_{\alpha,q}\in \mathbb{N}$ such that \[
\|\psi_{\alpha}\mathbf{F}-\mathbf{G}_{\alpha,k_{\alpha,q}}\|_{W^{m,p}(U_{\alpha},\mathbf{E}\restriction_{U_{\alpha}})}<\frac{1}{2^{\alpha+q}}.
\]

Since $\supp(\mathbf{G}_{\alpha,k})\subseteq U_{\alpha}$ for $\alpha,k\in\mathbb{N}$ and the family $(U_{\alpha})_{\alpha\in\mathbb{N}}$ is locally finite, the sum \[
\mathbf{G}_{q}:=\displaystyle\sum_{\alpha=1}^{\infty}\mathbf{G}_{\alpha,k_{\alpha,q}}
\] is locally finite and therefore defines a global smooth section $\mathbf{G}_{q}\in \Gamma(\mathbf{E})$.

Let us now show that $\mathbf{G}_q$ approximates $\mathbf{F}$ in the Sobolev norm. Using the triangle inequality gives \[
\|\mathbf{F}-\mathbf{G}_{q}\|_{W^{m,p}(M,\mathbf{E})}
\leq \displaystyle\sum_{\alpha=1}^{\infty}\|\psi_{\alpha}\mathbf{F}-\mathbf{G}_{\alpha,k_{\alpha,q}}\|_{W^{m,p}(M,\mathbf{E})}
\leq \displaystyle\sum_{\alpha=1}^{\infty}\frac{1}{2^{\alpha+q}}=\frac{1}{2^{q}}.
\] In particular, \[
\|\mathbf{G}_q\|_{W^{m,p}(M,\mathbf{E})}
\leq \|\mathbf{F}\|_{W^{m,p}(M,\mathbf{E})} + \|\mathbf{F}-\mathbf{G}_q\|_{W^{m,p}(M,\mathbf{E})}
\leq\|\mathbf{F}\|_{W^{m,p}(M,\mathbf{E})} + 2^{-q}<\infty,
\] so that $\mathbf{G}_q\in W^{m,p}(M,\mathbf{E})\cap \Gamma(\mathbf{E})$ for every $q\in \mathbb{N}$.

Since $2^{-q}\to 0$ as $q\to\infty$, we conclude that $\mathbf{G}_{q}\to \mathbf{F}$ in $W^{m,p}(M,\mathbf{E})$. That is, every section in $W^{m,p}(M,\mathbf{E})$ can be approximated in the Sobolev norm by smooth sections, showing that \[
H^{m,p}(M,\mathbf{E})=W^{m,p}(M,\mathbf{E}).
\] \end{proof}

\begin{theorem}[Chain rule for Lipschitz functions] \label{teo:regla-cadena-lipschitz-sobolev} \index{chain rule!for Sobolev functions} Let $(M,\mathbf{g})$ be a Riemannian manifold without boundary, let $1\leq p<\infty$, and let $u\in W^{1,p}_{\mathrm{loc}}(M)$ be a real-valued function. If $\Phi\colon \mathbb R\longrightarrow\mathbb R$ is $L$-Lipschitz, then $\Phi\circ u\in W^{1,p}_{\mathrm{loc}}(M)$ and \[
|d(\Phi\circ u)|_{\mathbf{g}}\leq L|du|_{\mathbf{g}}
\] almost everywhere. On the set of points $x$ for which $\Phi$ is differentiable at $u(x)$, we have \[
d(\Phi\circ u)_x=\Phi'(u(x))\,du_x
\] almost everywhere. Moreover, $du=0$ almost everywhere on each level set $\{u=c\}$.

If $u\in W^{1,p}(M)$ and $\Phi\circ u\in L^p(M)$, then $\Phi\circ u\in W^{1,p}(M)$. This integrability condition holds, in particular, when $\Phi(0)=0$ or when $\lambda_{\mathbf{g}}(M)<\infty$.

For $u,v\in W^{1,p}_{\mathrm{loc}}(M)$ we obtain the formulas \[
d(u-c)_+=\mathbf 1_{\{u>c\}}\,du,
\qquad
du_-=-\mathbf 1_{\{u<0\}}\,du,
\] and \[
d(u-v)_+
=
\mathbf 1_{\{u>v\}}(du-dv)
\] almost everywhere, where $t_+:=\displaystyle\max\{t,0\}$ and $u_-:=(-u)_+$. \end{theorem}

\begin{proof} We begin with a function $\Psi\in C^1(\mathbb R)$ whose derivative is bounded. Let $U\Subset M$ be a regular coordinate ball. Since $u\in W^{1,p}_{\mathrm{loc}}(M)$, the restriction $u|_U$ belongs to $W^{1,p}(U)$. Apply the Meyers--Serrin theorem for vector bundles, Theorem~\ref{meyers-serrin-haz}, to the trivial real bundle over $U$. This gives a sequence $(u_j)\subseteq C^\infty(U)\cap W^{1,p}(U)$ such that $u_j\to u$ in $W^{1,p}(U)$. After passing to a subsequence, we may also assume that $u_j\to u$ almost everywhere on $U$. Since $\Psi$ is Lipschitz, $\Psi(u_j)\to\Psi(u)$ in $L^p(U)$. For the differentials we have \[
d(\Psi\circ u_j)=\Psi'(u_j)\,du_j.
\] The inequality \[
\begin{split}
\|\Psi'(u_j)du_j-\Psi'(u)du\|_{L^p(U,T^*U)}
\leq{}&
\|\Psi'\|_{L^\infty(\mathbb R)}\|du_j-du\|_{L^p(U,T^*U)}\\
&+
\|(\Psi'(u_j)-\Psi'(u))du\|_{L^p(U,T^*U)}
\end{split}
\] shows that the left-hand side tends to zero. Indeed, the first term converges to zero by approximation in $W^{1,p}(U)$, and the second by the dominated convergence theorem~\ref{convergencia dominada}, since $\Psi'(u_j)\to\Psi'(u)$ almost everywhere and $|\Psi'(u_j)-\Psi'(u)|^p|du|_{\mathbf{g}}^p\leq 2^p\|\Psi'\|_{L^\infty(\mathbb R)}^p|du|_{\mathbf{g}}^p$. Closedness of the weak differential, Lemma~\ref{derivada debil es operador cerrado}, then implies that \[
d(\Psi\circ u)=\Psi'(u)\,du
\] on $U$. Since the regular coordinate balls with compact closure cover $M$, the identity holds locally throughout the manifold.

Before treating an arbitrary Lipschitz function, we prove the assertion about level sets. Fix $c\in\mathbb R$ and choose $\eta\in C_c^\infty(\mathbb R)$ such that $0\leq\eta\leq1$, $\eta=1$ on $[-1,1]$, and $\eta=0$ outside $[-2,2]$. For $\varepsilon>0$ define \[
\Psi_\varepsilon(t)
:=
\int_c^t\eta\left(\frac{s-c}{\varepsilon}\right)ds.
\] We have $|\Psi_\varepsilon|\leq2\varepsilon$ and $\Psi_\varepsilon'(t)=\eta\!\left(\frac{t-c}{\varepsilon}\right)$. By what we have already proved, \[
d(\Psi_\varepsilon\circ u)
=
\eta\left(\frac{u-c}{\varepsilon}\right)du.
\] Fix a regular coordinate ball $U\Subset M$. Since $|\Psi_\varepsilon\circ u|\leq2\varepsilon$, we have $\Psi_\varepsilon\circ u\to0$ in $L^p(U)$. On the other hand, \[
\eta\left(\frac{u-c}{\varepsilon}\right)du
\longrightarrow
\mathbf 1_{\{u=c\}}\,du
\] in $L^p(U,T^*U)$ by the dominated convergence theorem~\ref{convergencia dominada}. Indeed, the scalar factor converges pointwise to $\mathbf 1_{\{u=c\}}$ and its absolute value is bounded by one. Apply Lemma~\ref{derivada debil es operador cerrado} to any sequence $\varepsilon_j\to0^+$. The functions converge to zero and their differentials converge to $\mathbf 1_{\{u=c\}}du$, so the latter must be the weak differential of the zero function. Consequently, $\mathbf 1_{\{u=c\}}du=0$ almost everywhere on $U$. Since these balls cover $M$, the assertion holds throughout the manifold.

Now let $\Phi$ be an $L$-Lipschitz function. Take the one-dimensional mollifier from Theorem~\ref{teo:molificadores-euclidianos} and set $\Phi_j:=\Phi*\rho_{\frac{1}{j}}$. Then $\Phi_j\in C^\infty(\mathbb R)$, \[
\|\Phi_j-\Phi\|_{L^\infty(\mathbb R)}\longrightarrow0,
\qquad
|\Phi_j'|\leq L.
\] In particular, $\Phi_j(u)\to\Phi(u)$ in $L^p(Q)$ for every $Q\Subset M$. By the first part of the proof, \[
d(\Phi_j\circ u)=\Phi_j'(u)\,du.
\]

Again fix a regular coordinate ball $U\Subset M$ and write $a_j:=\Phi_j'(u)$. The sequence $(a_j)$ lies in the closed ball of radius $L$ in $L^\infty(U)$. Since $L^1(U)$ is separable by Lemma~\ref{Lp separable espacio de medida separable}, Corollary~\ref{cor: subsucesion debil estrella} allows us to pass to a subsequence, without changing notation, such that $a_j\rightharpoonup^*a$ in $L^\infty(U)$ for some function $a$ with $|a|\leq L$ almost everywhere. Let $\boldsymbol{\omega}\in\Gamma_c(T^*U)$. Since $\langle du,\boldsymbol{\omega}\rangle_{\mathbf{g}}\in L^1(U)$, weak-$*$ convergence gives \[
\int_Ua_j\langle du,\boldsymbol{\omega}\rangle_{\mathbf{g}}\,d\lambda_{\mathbf{g}}
\longrightarrow
\int_Ua\langle du,\boldsymbol{\omega}\rangle_{\mathbf{g}}\,d\lambda_{\mathbf{g}}.
\] For each $j$, the defining identity for the weak differential reads \[
\int_Ua_j\langle du,\boldsymbol{\omega}\rangle_{\mathbf{g}}\,d\lambda_{\mathbf{g}}
=
\int_U\Phi_j(u)d^*\boldsymbol{\omega}\,d\lambda_{\mathbf{g}}.
\] The right-hand side converges to $\displaystyle\int_U\Phi(u)d^*\boldsymbol{\omega}\,d\lambda_{\mathbf{g}}$, since $\Phi_j(u)\to\Phi(u)$ in $L^p(U)$ and $d^*\boldsymbol{\omega}$ is smooth with compact support. Therefore, \[
d(\Phi\circ u)=a\,du
\] on $U$. In particular, $|d(\Phi\circ u)|_{\mathbf{g}}\leq L|du|_{\mathbf{g}}$ almost everywhere. On intersections of two balls, the differentials thus obtained agree by uniqueness of the weak derivative. They therefore define the weak differential of $\Phi\circ u$ throughout $M$.

It remains to identify the coefficient $a$ at values where $\Phi$ is differentiable. If $\Phi$ is differentiable at $t$, then $\Phi_j'(t)\to\Phi'(t)$. To see this directly, write \[
\Phi_j'(t)
=
j\int_{\mathbb R}\Phi\left(t-\frac{s}{j}\right)\rho'(s)\,ds.
\] Write $\Phi(t+h)=\Phi(t)+h\Phi'(t)+r(h)$. For each $\varepsilon>0$ there is $\delta>0$ such that $|r(h)|\leq\varepsilon|h|$ if $|h|<\delta$. Moreover, \[
\Phi\left(t-\frac{s}{j}\right)
=
\Phi(t)-\frac{s}{j}\Phi'(t)+r\left(-\frac{s}{j}\right).
\] The remainder, multiplied by $j$, converges uniformly to zero for $s\in\operatorname{supp}\rho$. Indeed, if this support is contained in $[-R,R]$ and $j>\frac{R}{\delta}$, then \[
\sup_{s\in[-R,R]}j\left|r\left(-\frac{s}{j}\right)\right|
\leq
R\varepsilon.
\] Since $\displaystyle\int_{\mathbb R}\rho'(s)\,ds=0$ and $\displaystyle\int_{\mathbb R}s\rho'(s)\,ds=-1$, the preceding formula converges to $\Phi'(t)$. Let $D\subseteq U$ be the measurable set of points $x$ for which $\Phi$ is differentiable at $u(x)$. Then $a_j(x)\to\Phi'(u(x))$ on $D$. For each $\boldsymbol{\omega}\in\Gamma_c(T^*U)$, the function $\mathbf 1_D\langle du,\boldsymbol{\omega}\rangle_{\mathbf{g}}$ belongs to $L^1(U)$. Weak-$*$ convergence and, on the other hand, the dominated convergence theorem~\ref{convergencia dominada} on $D$ give \[
\int_Da\langle du,\boldsymbol{\omega}\rangle_{\mathbf{g}}\,d\lambda_{\mathbf{g}}
=
\int_D\Phi'(u)\langle du,\boldsymbol{\omega}\rangle_{\mathbf{g}}\,d\lambda_{\mathbf{g}}.
\] Since $\boldsymbol{\omega}$ is arbitrary, the localization principle implies that \[
a\,du=\Phi'(u)\,du
\] almost everywhere there. This proves the stated chain rule formula.

If $u\in W^{1,p}(M)$ and $\Phi(u)\in L^p(M)$, the bound on the differential shows that $d(\Phi(u))\in L^p(M,T^*M)$, so that $\Phi(u)\in W^{1,p}(M)$. If $\Phi(0)=0$, then $|\Phi(u)|\leq L|u|$. If $M$ has finite volume, use $|\Phi(u)|\leq|\Phi(0)|+L|u|$. In both cases this gives the required integrability.

Finally, the function $t\mapsto(t-c)_+$ is differentiable outside $c$, with derivative equal to $\mathbf 1_{(c,\infty)}$. The formula already proved holds outside $\{u=c\}$, and on this set $du=0$. Therefore, $d(u-c)_+=\mathbf 1_{\{u>c\}}du$. The identity for $u_-$ follows by applying this formula to $-u$, and the last identity follows by applying it to $u-v$. \end{proof}

\section{Sobolev embeddings and the Rellich--Kondrashov theorem for vector bundles} Sobolev embeddings are a basic tool for studying partial differential equations on bundles. On a compact manifold, a finite cover by regular coordinate balls reduces the problem to Euclidean embeddings for the components, after which the sections can be reconstructed using a partition of unity. Lemma~\ref{lema:meyers-serrin-haz-E} ensures that this procedure preserves Sobolev regularity exactly.

\begin{theorem}[Sobolev embedding for vector bundles] \label{teo: encaje de sobolev haces} \index{Sobolev embedding theorem@Sobolev embedding theorem!for vector bundles} Let $(M,\mathbf{g})$ be a closed Riemannian manifold of dimension $n$, and let $\mathbf{E}\longrightarrow M$ be a smooth real or complex vector bundle of finite rank, equipped with a bundle metric $\mathbf h_{\mathbf E}$ (Hermitian in the complex case) and a compatible connection $\nabla^{\mathbf E}$. Let $j\geq0$ and $m\geq1$ be integers and $1\leq p,q<\infty$ such that \[
j+m-\frac np\geq j-\frac nq.
\] Then the inclusion \[
W^{j+m,p}(M,\mathbf{E})\hookrightarrow W^{j,q}(M,\mathbf{E})
\] is continuous. \end{theorem}

\begin{proof} Fix the orders and exponents in the statement. Choose a finite cover $(U_\alpha,\phi_\alpha)_{\alpha=1}^N$ by regular coordinate balls and, on each $U_\alpha$, a frame $(\mathbf{e}_{1,\alpha},\dots,\mathbf{e}_{r,\alpha})$ defined on a neighborhood of $\overline{U_\alpha}$. Take a partition of unity $(\psi_\alpha)_{\alpha=1}^N$ subordinate to this cover. Since $M$ is compact, we may and do assume that \[
K_\alpha:=\operatorname{supp}(\psi_\alpha)\Subset U_\alpha
\qquad
\text{for every }\alpha\in\{1,\dots,N\}.
\] Set $\Omega_\alpha:=\phi_\alpha(U_\alpha)$.

Compatibility of $\nabla^{\mathbf E}$ with $\mathbf h_{\mathbf E}$ allows us to apply Lemma~\ref{lema:meyers-serrin-haz-E} and its local comparison estimates. For complex components, we apply the Euclidean estimates to their real and imaginary parts. Theorem~\ref{sobolev generalizado}, with total orders $j+m$ and $j$, gives, for each $\alpha$, a constant $C_\alpha>0$ such that \begin{equation}
\label{eq: contencion coordenadas para encaje de Sobolev}
\|v\|_{W^{j,q}(\Omega_\alpha)}
\leq
C_\alpha\|v\|_{W^{j+m,p}(\Omega_\alpha)}
\qquad
\text{for every }v\in W_0^{j+m,p}(\Omega_\alpha).
\end{equation} The formulation with $W_0^{j+m,p}(\Omega_\alpha)$ suffices: all functions appearing below have compact support in $\Omega_\alpha$.

Let $\mathbf{F}\in W^{j+m,p}(M,\mathbf{E})$ and write \[
\mathbf{F}=\sum_{a=1}^r F_\alpha^a \mathbf{e}_{a,\alpha}
\qquad\text{in }U_\alpha.
\] The localized section $\psi_\alpha \mathbf{F}$ has support contained in $K_\alpha$. By Lemma~\ref{lema:meyers-serrin-haz-E}, the functions \[
v_\alpha^a
:=
(\psi_\alpha F_\alpha^a)\circ\phi_\alpha^{-1}
\] belong to $W^{j+m,p}(\Omega_\alpha)$. Moreover, \[
\operatorname{supp}(v_\alpha^a)
\subseteq
\phi_\alpha(K_\alpha)
\Subset
\Omega_\alpha,
\] so that $v_\alpha^a\in W_0^{j+m,p}(\Omega_\alpha)$. Estimate \eqref{eq: contencion coordenadas para encaje de Sobolev} and the local comparison of norms give \begin{equation}
\label{eq:estimacion-local-encaje-haces-cap9}
\|(\psi_\alpha \mathbf{F})\restriction_{U_\alpha}\|_{W^{j,q}(U_\alpha,\mathbf{E})}
\leq
C_\alpha'
\|(\psi_\alpha \mathbf{F})\restriction_{U_\alpha}\|_{W^{j+m,p}(U_\alpha,\mathbf{E})}
\end{equation} for a constant $C_\alpha'>0$ independent of $\mathbf{F}$.

Let us make the global reconstruction of each localized section precise. For $s\in\{0,\dots,j\}$, let \[
G_{\alpha,s}
:=
\nabla_w^s\bigl((\psi_\alpha \mathbf{F})\restriction_{U_\alpha}\bigr).
\] Locality of the weak covariant derivative implies the inclusion \[
\operatorname{supp}_{U_\alpha}(G_{\alpha,s})
\subseteq
\operatorname{supp}_{U_\alpha}
\bigl((\psi_\alpha \mathbf{F})\restriction_{U_\alpha}\bigr)
\subseteq K_\alpha;
\] equality of supports is neither needed nor valid in general. Denote by $\widetilde G_{\alpha,s}$ the extension of $G_{\alpha,s}$ by zero to $M$. Since $K_\alpha\Subset U_\alpha$, this extension belongs to $L^q(M,T^{(0,s)}(TM)\otimes \mathbf{E})$.

Let us show that $\widetilde G_{\alpha,s}$ is the $s$th weak covariant derivative of $\psi_\alpha \mathbf{F}$. Choose $\eta_\alpha\in C_c^\infty(U_\alpha)$ equal to one on a neighborhood of $K_\alpha$. For $\boldsymbol{\Phi}\in\Gamma_c^\infty(T^{(0,s)}(TM)\otimes \mathbf{E})$, the section $\eta_\alpha\boldsymbol{\Phi}\restriction_{U_\alpha}$ is a test function on $U_\alpha$. The definition of $G_{\alpha,s}$ and locality of $(\nabla^s)^*$ give \[
\begin{aligned}
\int_M
\langle\widetilde G_{\alpha,s},\boldsymbol{\Phi}\rangle_{\mathbf{g},\mathbf{h}_{\mathbf{E}}}\,d\lambda_{\mathbf{g}}
&=
\int_{U_\alpha}
\langle G_{\alpha,s},\eta_\alpha\boldsymbol{\Phi}\rangle_{\mathbf{g},\mathbf{h}_{\mathbf{E}}}\,d\lambda_{\mathbf{g}}\\
&=
\int_{U_\alpha}
\left\langle
\psi_\alpha \mathbf{F},(\nabla^s)^*(\eta_\alpha\boldsymbol{\Phi})
\right\rangle_{\mathbf{h}_{\mathbf{E}}}\,d\lambda_{\mathbf{g}}\\
&=
\int_M
\langle\psi_\alpha \mathbf{F},(\nabla^s)^*\boldsymbol{\Phi}\rangle_{\mathbf{h}_{\mathbf{E}}}\,d\lambda_{\mathbf{g}}.
\end{aligned}
\] In the last equality we used $\eta_\alpha=1$ on a neighborhood of the support of $\psi_\alpha \mathbf{F}$. Therefore, $\psi_\alpha \mathbf{F}\in W^{j,q}(M,\mathbf{E})$ and its global norm agrees with the norm obtained by integrating over $U_\alpha$.

Finally, we need a bound for localization in the source space. For smooth sections, the iterated Leibniz rule gives, for $0\leq s\leq j+m$, \[
|\nabla^s(\psi_\alpha \mathbf{F})|_{\mathbf{g},\mathbf{h}_{\mathbf{E}}}
\leq
A_s\sum_{t=0}^s
|\nabla^{s-t}\psi_\alpha|_{\mathbf{g}}
|\nabla^t\mathbf{F}|_{\mathbf{g},\mathbf{h}_{\mathbf{E}}}.
\] The derivatives of $\psi_\alpha$ are bounded because $M$ is compact. Summing the $L^p(\Omega_\alpha)$ norms gives a constant $A_\alpha>0$ such that \[
\|\psi_\alpha \mathbf{F}\|_{W^{j+m,p}(M,\mathbf{E})}
\leq
A_\alpha\|\mathbf{F}\|_{W^{j+m,p}(M,\mathbf{E})}
\] for every smooth section $\mathbf{F}$. The Meyers--Serrin theorem and passage to the limit in distributions extend the inequality to every $\mathbf{F}\in W^{j+m,p}(M,\mathbf{E})$. Since the cover is finite, there is $A>0$ such that \begin{equation}
\label{eq: acotacion particion de la unidad y norma sobolev auxiliar para encaje de Sobolev y Rellich--Kondrashov}
\|\psi_\alpha \mathbf{F}\|_{W^{j+m,p}(M,\mathbf{E})}
\leq
A\|\mathbf{F}\|_{W^{j+m,p}(M,\mathbf{E})}
\end{equation} for every $\alpha\in\{1,\dots,N\}$.

Using $\mathbf{F}=\displaystyle\sum_{\alpha=1}^N\psi_\alpha \mathbf{F}$, the triangle inequality, \eqref{eq:estimacion-local-encaje-haces-cap9}, and \eqref{eq: acotacion particion de la unidad y norma sobolev auxiliar para encaje de Sobolev y Rellich--Kondrashov} gives \[
\begin{aligned}
\|\mathbf{F}\|_{W^{j,q}(M,\mathbf{E})}
&\leq
\sum_{\alpha=1}^N
\|\psi_\alpha \mathbf{F}\|_{W^{j,q}(M,\mathbf{E})}\\
&\leq
\sum_{\alpha=1}^N
C_\alpha'
\|\psi_\alpha \mathbf{F}\|_{W^{j+m,p}(M,\mathbf{E})}\\
&\leq
C\|\mathbf{F}\|_{W^{j+m,p}(M,\mathbf{E})}.
\end{aligned}
\] This proves continuity of the embedding. \end{proof}

\begin{theorem}[Rellich--Kondrashov in vector bundles] \label{teo: rellich--kondrashov para haces vectoriales} \index{Rellich Kondrashov in vector bundles@Rellich--Kondrashov in vector bundles} Let $(M,\mathbf{g})$ be a compact Riemannian manifold without boundary of dimension $n$, and let $\mathbf{E}\to M$ be a smooth real or complex vector bundle of finite rank, equipped with a bundle metric $\mathbf{h}_{\mathbf{E}}$ (Hermitian in the complex case) and a compatible connection $\nabla^{\mathbf{E}}$. Let $j\in\mathbb N_0$, let $m\geq1$ be an integer, and let $1\leq p,q<\infty$ satisfy \[
j+m-\frac{n}{p}>j-\frac{n}{q}.
\] Then the inclusion \[
W^{j+m,p}(M,\mathbf{E})\hookrightarrow W^{j,q}(M,\mathbf{E})
\] is compact. \end{theorem}

\begin{proof} Let $(\mathbf{F}_\nu)_{\nu\in\mathbb N}$ be a bounded sequence in $W^{j+m,p}(M,\mathbf{E})$. Compatibility of $\nabla^{\mathbf E}$ with $\mathbf h_{\mathbf E}$ allows us to apply the local comparison of norms in Lemma~\ref{lema:meyers-serrin-haz-E}. In the complex case, consider the underlying real bundle with metric $\operatorname{Re}\mathbf h_{\mathbf E}$ and the induced real connection; its pointwise norm agrees with the Hermitian norm, and its components are the real and imaginary parts of the complex components. Retain the cover, frames, and functions $\psi_\alpha$ from the preceding proof. By \eqref{eq: acotacion particion de la unidad y norma sobolev auxiliar para encaje de Sobolev y Rellich--Kondrashov}, for each $\alpha$ the sequence $(\psi_\alpha \mathbf{F}_\nu)_\nu$ is bounded in $W^{j+m,p}(M,\mathbf{E})$. Lemma~\ref{lema:meyers-serrin-haz-E} then shows that, for each pair $(\alpha,a)$, the sequence \[
v_{\alpha,\nu}^a
:=
\bigl(\psi_\alpha(\mathbf{F}_\nu)_\alpha^a\bigr)
\circ\phi_\alpha^{-1}
\] is bounded in $W^{j+m,p}(\Omega_\alpha)$. Moreover, \[
\operatorname{supp}(v_{\alpha,\nu}^a)
\subseteq
\phi_\alpha(K_\alpha)
\Subset
\Omega_\alpha
\] for every $\nu$. In particular, $v_{\alpha,\nu}^a\in W_0^{j+m,p}(\Omega_\alpha)$.

For each pair $(\alpha,a)$, Theorem~\ref{rellich kondrashov generalizado} gives a subsequence converging in $W^{j,q}(\Omega_\alpha)$. Since there are only $Nr$ pairs, a finite diagonal extraction produces a subsequence, which we do not relabel, such that $(v_{\alpha,\nu}^a)_\nu$ is Cauchy in $W^{j,q}(\Omega_\alpha)$ for every $\alpha$ and every $a$.

The local comparison of norms and reconstruction by extension by zero used in the proof of the preceding theorem give, for $\nu,\mu\in\mathbb N$, \[
\begin{aligned}
\|\mathbf{F}_\nu-\mathbf{F}_\mu\|_{W^{j,q}(M,\mathbf{E})}
&\leq
\sum_{\alpha=1}^N
\|\psi_\alpha(\mathbf{F}_\nu-\mathbf{F}_\mu)\|_{W^{j,q}(M,\mathbf{E})}\\
&\leq
C\sum_{\alpha=1}^N\sum_{a=1}^r
\|v_{\alpha,\nu}^a-v_{\alpha,\mu}^a\|_{W^{j,q}(\Omega_\alpha)}.
\end{aligned}
\] The right-hand side tends to zero as $\nu,\mu\to\infty$. The subsequence is therefore Cauchy in the Banach space $W^{j,q}(M,\mathbf{E})$ and converges in it. This proves compactness. \end{proof}

\begin{corollary}[Embeddings and the critical exponent] \label{cor:encajes-criticos-haces} Under the geometric hypotheses of the preceding theorems, let $j\geq0$ and $m\geq1$ be integers, and let $1\leq p<n/m$. Then \[
W^{j+m,p}(M,\mathbf E)\hookrightarrow W^{j,q}(M,\mathbf E)
\] continuously for $1\leq q\leq p_m^*:=\displaystyle\frac{np}{n-mp}$, and compactly for $1\leq q<p_m^*$. \end{corollary} \begin{proof} Since $mp<n$, the condition $j+m-\displaystyle\frac np\geq j-\displaystyle\frac nq$ is equivalent to $q\leq p_m^*$, and the strict inequality is equivalent to $q<p_m^*$. Apply Theorems~\ref{teo: encaje de sobolev haces} and \ref{teo: rellich--kondrashov para haces vectoriales}, respectively. \end{proof}

\subsection{Embeddings into spaces of Hölder sections}

Solutions of many variational problems are initially obtained as elements of a Sobolev space. The localization principle for vector bundles in Lemma~\ref{lema: principio de localizacion para haces vectoriales} allows us to convert the variational identity into a local equation in the weak sense. This step alone, however, does not yield classical regularity. The next question is when the weak solution admits a differentiable representative and with what continuity its derivatives vary. As in the Euclidean case, the embedding does not replace estimates specific to the equation: once those estimates give a bound in $W^{m,p}(M,\mathbf{E})$, spaces of Hölder sections translate this Sobolev regularity into quantitative classical regularity.

Before defining these spaces, we must specify how to compare vectors belonging to different fibers. We use the notation for parallel transport fixed in Definition~\ref{def:variedades-riemannianas-transporte-paralelo}. Throughout this subsection the manifolds have no boundary. Initially we assume that $M$ is connected, because the Riemannian distance $d_{\mathbf{g}}$ was defined under this hypothesis. If $M$ is disconnected, the definitions and results apply componentwise. In particular, connectedness imposes no restriction for compact manifolds, which have finitely many connected components.

\begin{definition}[Fiberwise separation by parallel transport] \label{def:separacion-fibrada-transporte-paralelo} \index{parallel transport!fiberwise separation} Let $(M,\mathbf{g})$ be a connected Riemannian manifold without boundary such that $\operatorname{inj}(M)>0$. Let $\mathbf{F}\longrightarrow M$ be a smooth vector bundle equipped with a bundle metric $\mathbf{h}_{\mathbf{F}}$ and a compatible connection $\nabla^{\mathbf{F}}$. Fix \[
\rho_M:=\min\{1,\operatorname{inj}(M)\}.
\] If $x,y\in M$ satisfy $0<d_{\mathbf{g}}(x,y)<\rho_M$, denote by $\gamma_{xy}\colon [0,1]\longrightarrow M$ the unique minimizing geodesic parametrized with constant speed satisfying \[
\gamma_{xy}(0)=x,
\qquad
\gamma_{xy}(1)=y.
\] For $\xi\in \mathbf{F}_x$ and $\eta\in \mathbf{F}_y$, define \begin{equation}
\label{eq:separacion-fibrada-transporte}
|\xi-\eta|_{\mathbf{h}_{\mathbf{F}};x,y}
:=
\left|
\xi-P^{\mathbf{F}}_{\gamma_{xy},1\to0}(\eta)
\right|_{\mathbf{h}_{\mathbf{F}}(x)}.
\end{equation} Since $\nabla^{\mathbf{F}}$ is compatible with $\mathbf{h}_{\mathbf{F}}$, parallel transport is an isometry. Proposition~\ref{prop: transporte paralelo isomorfismo} then gives the equivalent expression \[
|\xi-\eta|_{\mathbf{h}_{\mathbf{F}};x,y}
=
\left|
P^{\mathbf{F}}_{\gamma_{xy},0\to1}(\xi)-\eta
\right|_{\mathbf{h}_{\mathbf{F}}(y)}.
\] When $x=y$, set \[
|\xi-\eta|_{\mathbf{h}_{\mathbf{F}};x,x}:=|\xi-\eta|_{\mathbf{h}_{\mathbf{F}}(x)}.
\]

If $\mathbf{T}\in\Gamma^0(\mathbf{F})$ and $S\subseteq M$ satisfies $\operatorname{diam}_{\mathbf{g}}(S)<\rho_M$, define the fiberwise oscillation of $\mathbf{T}$ on $S$ by \[
\operatorname{osc}_{\mathbf{F}}(\mathbf{T};S)
:=
\sup_{x,y\in S}|\mathbf{T}(x)-\mathbf{T}(y)|_{\mathbf{h}_{\mathbf{F}};x,y}.
\] For $0<\alpha<1$, define \begin{equation}
\label{eq:seminorma-holder-fibrada}
[\mathbf{T}]_{\alpha,\mathbf{F}}
:=
\sup_{\substack{x,y\in M\\0<d_{\mathbf{g}}(x,y)<\rho_M}}
\frac{|\mathbf{T}(x)-\mathbf{T}(y)|_{\mathbf{h}_{\mathbf{F}};x,y}}{d_{\mathbf{g}}(x,y)^\alpha}.
\end{equation} If $U\subseteq M$, write $[\mathbf{T}]_{\alpha,\mathbf{F};U}$ when the supremum is taken only over pairs $x,y\in U$. \end{definition}

\begin{remark}[Choice of radius] \label{obs:eleccion-radio-holder-secciones} The bound by $1$ in the definition of $\rho_M$ merely normalizes the scale. If $0<\rho_1\leq\rho_2\leq\operatorname{inj}(M)$ and $[\mathbf{T}]_{\alpha,\mathbf{F};\rho_j}$ denotes the seminorm obtained using radius $\rho_j$, then \[
[\mathbf{T}]_{\alpha,\mathbf{F};\rho_1}
\leq
[\mathbf{T}]_{\alpha,\mathbf{F};\rho_2}
\] and \[
[\mathbf{T}]_{\alpha,\mathbf{F};\rho_2}
\leq
[\mathbf{T}]_{\alpha,\mathbf{F};\rho_1}
+
2\rho_1^{-\alpha}\|\mathbf{T}\|_{L^\infty(M,\mathbf{F})}.
\] The second inequality follows by separating pairs at distance less than $\rho_1$ from those whose distance belongs to $[\rho_1,\rho_2)$. Thus, after adding the uniform norm, any two positive radii less than or equal to $\operatorname{inj}(M)$ produce equivalent norms.

One may also use the formulation in terms of balls \[
[\mathbf{T}]_{\alpha,\mathbf{F}}^{\mathrm{bol}}
:=
\sup_{\substack{x\in M\\0<r<\frac{\rho_M}{2}}}
\frac{\operatorname{osc}_{\mathbf{F}}(\mathbf{T};B_{\mathbf{g}}(x,r))}{r^\alpha}.
\] If $x,y\in B_{\mathbf{g}}(z,r)$, then $d_{\mathbf{g}}(x,y)<2r$, and therefore \[
[\mathbf{T}]_{\alpha,\mathbf{F}}^{\mathrm{bol}}
\leq
2^\alpha[\mathbf{T}]_{\alpha,\mathbf{F}}.
\] Conversely, if $d_{\mathbf{g}}(x,y)<\frac{\rho_M}{2}$, take $d_{\mathbf{g}}(x,y)<r<\frac{\rho_M}{2}$ and use $y\in B_{\mathbf{g}}(x,r)$. Letting $r\to d_{\mathbf{g}}(x,y)^{+}$ gives \[
\frac{|\mathbf{T}(x)-\mathbf{T}(y)|_{\mathbf{h}_{\mathbf{F}};x,y}}{d_{\mathbf{g}}(x,y)^\alpha}
\leq
[\mathbf{T}]_{\alpha,\mathbf{F}}^{\mathrm{bol}}.
\] For $\frac{\rho_M}{2}\leq d_{\mathbf{g}}(x,y)<\rho_M$, we have \[
\frac{|\mathbf{T}(x)-\mathbf{T}(y)|_{\mathbf{h}_{\mathbf{F}};x,y}}{d_{\mathbf{g}}(x,y)^\alpha}
\leq
2^{\alpha+1}\rho_M^{-\alpha}\|\mathbf{T}\|_{L^\infty(M,\mathbf{F})}.
\] Thus the quantities \[
\|\mathbf{T}\|_{L^\infty(M,\mathbf{F})}+[\mathbf{T}]_{\alpha,\mathbf{F}}
\qquad\text{and}\qquad
\|\mathbf{T}\|_{L^\infty(M,\mathbf{F})}+[\mathbf{T}]_{\alpha,\mathbf{F}}^{\mathrm{bol}}
\] are equivalent. The factor $\frac{1}{2}$ ensures that any two points of $B_{\mathbf{g}}(x,r)$ lie at distance less than $\rho_M$. Likewise, uniqueness of the minimizing geodesic requires $d_{\mathbf{g}}(x,y)<\operatorname{inj}(M)$; the non-strict inequality is insufficient. \end{remark}

\begin{definition}[Hölder sections and Hölder sections up to the closure] \label{def:secciones-holder-haz} \index{Holder section@Hölder section} Let $(M,\mathbf{g})$ be a connected Riemannian manifold without boundary with $\operatorname{inj}(M)>0$, and let $\mathbf{E}\longrightarrow M$ be a smooth vector bundle equipped with a bundle metric $\mathbf{h}_{\mathbf{E}}$ and a compatible connection $\nabla^{\mathbf{E}}$. For $j\in\mathbb N_0$ set \[
\mathbf{E}_j:=T^{(0,j)}(TM)\otimes \mathbf{E},
\] with the convention $\mathbf{E}_0=\mathbf{E}$. On $\mathbf{E}_j$ consider the bundle metric induced by $\mathbf{g}$ and $\mathbf{h}_{\mathbf{E}}$, and the product connection induced by the Levi--Civita connection and $\nabla^{\mathbf{E}}$.

Let $k\in\mathbb N_0$ and $0<\alpha<1$. For $\mathbf{u}\in\Gamma^k(\mathbf{E})$ define \[
\|\nabla^j \mathbf{u}\|_{L^\infty(M,\mathbf{E}_j)}
:=
\sup_{x\in M}|\nabla^j \mathbf{u}(x)|_{\mathbf{g},\mathbf{h}_{\mathbf{E}}},
\qquad
0\leq j\leq k,
\] and \begin{equation}
\label{eq:norma-gamma-k-alpha}
\|\mathbf{u}\|_{\Gamma^{k,\alpha}(\mathbf{E})}
:=
\sum_{j=0}^{k}\|\nabla^j \mathbf{u}\|_{L^\infty(M,\mathbf{E}_j)}
+
[\nabla^k \mathbf{u}]_{\alpha,\mathbf{E}_k}.
\end{equation} The space of Hölder sections of order $(k,\alpha)$ is \[
\Gamma^{k,\alpha}(\mathbf{E})
:=
\left\{
 \mathbf{u}\in\Gamma^k(\mathbf{E})
 \middle|
 \|\mathbf{u}\|_{\Gamma^{k,\alpha}(\mathbf{E})}<\infty
\right\}.
\]

Now let $U\Subset M$ be open. For $\mathbf{u}\in\Gamma^k(\mathbf{E}\restriction_U)$ set \[
\mathcal E_{\overline U}^{k,\alpha}(\mathbf{u})
:=
\left\{
 \widetilde{\mathbf{u}}\in\Gamma^{k,\alpha}(\mathbf{E})
 \middle|
 \widetilde{\mathbf{u}}\restriction_U=\mathbf{u}
\right\}.
\] Define \[
\Gamma^{k,\alpha}(\overline U,\mathbf{E})
:=
\left\{
 \mathbf{u}\in\Gamma^k(\mathbf{E}\restriction_U)
 \middle|
 \mathcal E_{\overline U}^{k,\alpha}(\mathbf{u})\neq\varnothing
\right\}
\] and equip this space with the extension norm \begin{equation}
\label{eq:norma-gamma-holder-cerradura-extension}
\|\mathbf{u}\|_{\Gamma^{k,\alpha}(\overline U,\mathbf{E})}
:=
\inf_{\widetilde{\mathbf{u}}\in\mathcal E_{\overline U}^{k,\alpha}(\mathbf{u})}
\|\widetilde{\mathbf{u}}\|_{\Gamma^{k,\alpha}(\mathbf{E})}.
\end{equation} The second argument records the ambient bundle. We are not regarding $\overline U$ as a smooth manifold. Since $U$ is dense in $\overline U$, an extension uniquely determines the values of $\mathbf{u}$ and its covariant derivatives through order $k$ on the closure.

If $M$ has finitely many connected components $M_1,\dots,M_N$, apply the definition on each component and write \[
\Gamma^{k,\alpha}(\mathbf{E})
:=
\left\{
 \mathbf{u}\in\Gamma^k(\mathbf{E})
 \middle|
 \mathbf{u}\restriction_{M_\nu}\in
 \Gamma^{k,\alpha}(\mathbf{E}\restriction_{M_\nu})
 \text{ for every }\nu\in\{1,\dots,N\}
\right\},
\] with the norm \[
\|\mathbf{u}\|_{\Gamma^{k,\alpha}(\mathbf{E})}
:=
\max_{1\leq\nu\leq N}
\|\mathbf{u}\restriction_{M_\nu}\|_{\Gamma^{k,\alpha}(\mathbf{E}\restriction_{M_\nu})}.
\] If \[
\|\mathbf{u}\|_{\Sigma}:=
\sum_{\nu=1}^{N}
\|\mathbf{u}\restriction_{M_\nu}\|_{\Gamma^{k,\alpha}(\mathbf{E}\restriction_{M_\nu})},
\] then \[
\|\mathbf{u}\|_{\Gamma^{k,\alpha}(\mathbf{E})}
\leq\|\mathbf{u}\|_{\Sigma}
\leq N\|\mathbf{u}\|_{\Gamma^{k,\alpha}(\mathbf{E})}.
\] This is the convention used for disconnected compact manifolds. \end{definition}

\begin{remark}[Local extensions of sections] \label{obs:extensiones-locales-secciones-holder} The global extension condition in the preceding definition is equivalent to the existence of an extension to an open neighborhood of $\overline U$. More precisely, it suffices that there be open sets \[
\overline U\subseteq V_0\Subset V_1\Subset V
\] and a section $\widehat{\mathbf{u}}\in\Gamma^k(\mathbf{E}\restriction_V)$ extending $\mathbf{u}$ and satisfying \[
\sum_{j=0}^{k}
\|\nabla^j\widehat{\mathbf{u}}\|_{L^\infty(V_1,\mathbf{E}_j)}
+
[\nabla^k\widehat{\mathbf{u}}]_{\alpha,\mathbf{E}_k;V_1}
<\infty.
\] Indeed, choose a cutoff function $\chi\in C_c^\infty(V_1)$ equal to one on a neighborhood of $\overline{V_0}$. The extension of $\chi\widehat{\mathbf{u}}$ by zero outside $V_1$, denoted by $\widetilde{\chi\widehat{\mathbf{u}}}$, belongs to $\Gamma^{k,\alpha}(\mathbf{E})$ and extends $\mathbf{u}$. The presence of $V_1$ ensures that all terms in the Leibniz rule are controlled on a compact subset of the extension domain. This is the analogue for sections of the convention used in $C^{k,\alpha}(\overline\Omega)$.

The term $\|\nabla^k \mathbf{u}\|_{L^\infty(M,T^{(0,k)}(TM)\otimes \mathbf{E})}$ cannot be omitted from \eqref{eq:norma-gamma-k-alpha}. The Hölder seminorm does not control the size of $\nabla^k \mathbf{u}$ and vanishes, for example, on parallel tensors. \end{remark}

The local comparison between the geometric norm and the scalar components uses two facts: uniform equivalence of the metric with the Euclidean metric on a regular coordinate ball, and control of parallel transport along short geodesics.

\begin{lemma}[Local comparison of the Hölder seminorm] \label{lem:comparacion-local-holder-transporte-componentes} Let $(M,\mathbf{g})$ be a connected Riemannian manifold without boundary with $\operatorname{inj}(M)>0$. Let $(U,\phi)$ be a regular coordinate ball in $M$. Let $\mathbf{F}\longrightarrow M$ be a smooth vector bundle of rank $N$, equipped with a bundle metric $\mathbf{h}_{\mathbf{F}}$ and a compatible connection $\nabla^{\mathbf{F}}$, and let $(\mathbf{f}_1,\dots,\mathbf{f}_N)$ be a local frame of $\mathbf{F}$ defined on a neighborhood of $\overline U$. If \[
\mathbf{T}=\sum_{A=1}^{N}T^A\mathbf{f}_A
\qquad\text{on $U$},
\] then \[
\mathbf{T}\in\Gamma^{0,\alpha}(\overline U,\mathbf{F})
\quad\Longleftrightarrow\quad
T^A\circ\phi^{-1}
\in C^{0,\alpha}(\overline{\phi(U)})
\] for every $A\in\{1,\dots,N\}$. Moreover, there are constants $c_\alpha,C_\alpha>0$, independent of $\mathbf{T}$, such that \begin{equation}
\label{eq:equivalencia-holder-local-componentes}
\begin{aligned}
c_\alpha
\sum_{A=1}^{N}
\|T^A\circ\phi^{-1}\|_{C^{0,\alpha}(\overline{\phi(U)})}
&\leq
\|\mathbf{T}\|_{\Gamma^{0,\alpha}(\overline U,\mathbf{F})}\\
&\leq
C_\alpha
\sum_{A=1}^{N}
\|T^A\circ\phi^{-1}\|_{C^{0,\alpha}(\overline{\phi(U)})}.
\end{aligned}
\end{equation} \end{lemma}

\begin{proof} Let $U'$ be a larger regular coordinate ball, so that $\overline U\subseteq U'$ and the chart and frame are defined on $U'$. Choose an open set $V$ with \[
\overline U\subseteq V\Subset U'.
\] Since $\overline V$ is compact, there is $0<\delta<\rho_M$ such that every minimizing geodesic of length less than $\delta$ with endpoints in $\overline V$ remains in $U'$. Decreasing $\delta$ if necessary, the Euclidean segment joining $\phi(x)$ to $\phi(y)$ remains in $\phi(U')$ whenever $x,y\in\overline V$ and $d_{\mathbf{g}}(x,y)<\delta$.

For $x,y\in\overline V$ with $d_{\mathbf{g}}(x,y)<\delta$, define the matrix $G(x,y)$ by \begin{equation}
\label{eq:matriz-transporte-marco-holder}
P^{\mathbf{F}}_{\gamma_{xy},1\to0}\bigl(\mathbf{f}_B(y)\bigr)
=
\sum_{A=1}^{N}G^A_{\ B}(x,y)\mathbf{f}_A(x).
\end{equation} Let $R(t)$ be the matrix of the transport $P^{\mathbf{F}}_{\gamma_{xy},t\to0}$ in the frame. By Proposition~\ref{prop: ecuacion local transporte paralelo operador}, $R(0)=I$ and $R$ satisfies a linear system whose coefficients have norm bounded by $C|\dot\gamma_{xy}(t)|_{\mathbf{g}}=Cd_{\mathbf{g}}(x,y)$. Grönwall's lemma~\ref{lema: gronwall} gives \[
\|R(t)\|_{\operatorname{op}(\mathbb K^N)}\leq e^{Cd_{\mathbf{g}}(x,y)t}
\] and, integrating the equation, \begin{equation}
\label{eq:matriz-transporte-cerca-identidad}
\|G(x,y)-I\|_{\operatorname{op}(\mathbb K^N)}
=
\|R(1)-R(0)\|_{\operatorname{op}(\mathbb K^N)}
\leq
C d_{\mathbf{g}}(x,y).
\end{equation} Here and below the constant may increase, but is independent of $x$ and $y$.

Uniform equivalence of $\mathbf{g}$ and the Euclidean metric in the chart gives constants $c_0,C_0>0$ such that \begin{equation}
\label{eq:comparacion-distancias-carta-holder}
c_0|\phi(x)-\phi(y)|
\leq
d_{\mathbf{g}}(x,y)
\leq
C_0|\phi(x)-\phi(y)|
\end{equation} for the preceding pairs. The first inequality follows by estimating the length of the minimizing geodesic in coordinates. For the second, estimate the length of the curve whose coordinate expression is the segment joining $\phi(x)$ to $\phi(y)$.

For $z\in U'$, let \[
\mathcal F_z\colon \mathbb K^N\longrightarrow \mathbf{F}_z,
\qquad
\mathcal F_z(c^1,\dots,c^N):=\sum_{A=1}^{N}c^A\mathbf{f}_A(z).
\] If $\mathbf{S}=\displaystyle\sum_{A=1}^{N}S^A\mathbf{f}_A$ and $s(z):=(S^1(z),\dots,S^N(z))$, then \eqref{eq:matriz-transporte-marco-holder} implies \[
\mathbf{S}(x)-P^{\mathbf{F}}_{\gamma_{xy},1\to0}(\mathbf{S}(y))
=
\mathcal F_x\bigl(s(x)-G(x,y)s(y)\bigr).
\] Therefore, \[
s(x)-G(x,y)s(y)
=
s(x)-s(y)+(I-G(x,y))s(y).
\] Uniform equivalence of the Euclidean norm of the components and $\mathbf{h}_{\mathbf{F}}$, together with \eqref{eq:matriz-transporte-cerca-identidad} and \eqref{eq:comparacion-distancias-carta-holder}, gives \begin{equation}
\label{eq:comparacion-holder-intrinseca-componentes}
\begin{aligned}
&\|\mathbf{S}\|_{L^\infty(\overline U,\mathbf{F})}+[\mathbf{S}]_{\alpha,\mathbf{F};\overline U}\\
&\qquad\leq
C\sum_{A=1}^{N}
\left(
\|S^A\circ\phi^{-1}\|_{C^0(\overline{\phi(U)})}
+
[S^A\circ\phi^{-1}]_{C^{0,\alpha}(\overline{\phi(U)})}
\right).
\end{aligned}
\end{equation} For the reverse inequality, solve for \[
s(x)-s(y)
=
\bigl(s(x)-G(x,y)s(y)\bigr)
+
(G(x,y)-I)s(y)
\] and repeat the same estimates. Pairs with $d_{\mathbf{g}}(x,y)\geq\delta$ are controlled by the uniform norms. This proves equivalence of the intrinsic norms on $\overline U$.

We turn to the extension norms. Let $\widetilde{\mathbf{T}}\in\mathcal E_{\overline U}^{0,\alpha}(\mathbf{T})$. Choose a cutoff function $\chi\in C_c^\infty(U')$ equal to one on a neighborhood of $\overline U$. The components of $\chi\widetilde{\mathbf{T}}$ in the frame, transported by $\phi$ and extended by zero outside $\phi(U')$, are global extensions of $T^A\circ\phi^{-1}$. Estimate \eqref{eq:comparacion-holder-intrinseca-componentes}, applied on a compact set containing $\operatorname{supp}\chi$, and the product rule give \[
\sum_{A=1}^{N}
\|T^A\circ\phi^{-1}\|_{C^{0,\alpha}(\overline{\phi(U)})}
\leq
C\|\widetilde{\mathbf{T}}\|_{\Gamma^{0,\alpha}(\mathbf{F})}.
\] Taking the infimum over $\widetilde{\mathbf{T}}$ yields the first inequality in \eqref{eq:equivalencia-holder-local-componentes}.

Conversely, choose extensions $\widetilde t^A\in C^{0,\alpha}(\mathbb R^n)$ of $T^A\circ\phi^{-1}$. On $U'$ define \[
\widehat{\mathbf{T}}
:=
\chi\sum_{A=1}^{N}(\widetilde t^A\circ\phi)\mathbf{f}_A
\] and denote its extension by zero outside $U'$ by $\widetilde{\widehat{\mathbf{T}}}$. This section belongs to $\Gamma^{0,\alpha}(\mathbf{F})$, extends $\mathbf{T}$, and satisfies \[
\|\widetilde{\widehat{\mathbf{T}}}\|_{\Gamma^{0,\alpha}(\mathbf{F})}
\leq
C\sum_{A=1}^{N}
\|\widetilde t^A\|_{C^{0,\alpha}(\mathbb R^n)}.
\] Taking the infimum first over the scalar extensions and then in the definition of the norm of $\mathbf{T}$ gives the second inequality. The two constructions also prove equivalence of membership. \end{proof}

\begin{lemma}[Multiplication by smooth functions in spaces of Hölder sections] \label{lem:multiplicacion-funcion-suave-holder-secciones} Let $(M,\mathbf{g})$ be a connected Riemannian manifold without boundary with $\operatorname{inj}(M)>0$, and let $\mathbf{E}\longrightarrow M$ be a smooth vector bundle equipped with a bundle metric and a compatible connection. Let $k\in\mathbb N_0$, $0<\alpha<1$, and $\chi\in C_c^\infty(M)$. Then there is a constant $C_{\chi,k,\alpha}>0$ such that \[
\|\chi \mathbf{u}\|_{\Gamma^{k,\alpha}(\mathbf{E})}
\leq
C_{\chi,k,\alpha}\|\mathbf{u}\|_{\Gamma^{k,\alpha}(\mathbf{E})}
\] for every $\mathbf{u}\in\Gamma^{k,\alpha}(\mathbf{E})$. \end{lemma}

\begin{proof} We first record two elementary estimates. Let $\mathbf{F}\to M$ be a bundle with a compatible connection, and let $\mathbf{S}\in\Gamma^1(\mathbf{F})$. For $x,y\in M$ with $0<d_{\mathbf{g}}(x,y)<\rho_M$, Lemma~\ref{lema: TFC covariante}, applied to $\gamma_{xy}\colon [0,1]\longrightarrow M$, gives \[
P^{\mathbf{F}}_{\gamma_{xy},1\to0}(\mathbf{S}(y))-\mathbf{S}(x)
=
\int_0^1
P^{\mathbf{F}}_{\gamma_{xy},t\to0}
\bigl(\nabla^{\mathbf{F}}_{\dot\gamma_{xy}(t)}\mathbf{S}\bigr)\,dt.
\] Since $|\dot\gamma_{xy}(t)|_{\mathbf{g}}=d_{\mathbf{g}}(x,y)$ and parallel transport is an isometry, \begin{equation}
\label{eq:seccion-c1-es-lipschitz-paralela-multiplicacion}
|\mathbf{S}(x)-\mathbf{S}(y)|_{\mathbf{h}_{\mathbf{F}};x,y}
\leq
d_{\mathbf{g}}(x,y)
\bigl\||\nabla^{\mathbf{F}}\mathbf{S}|_{\mathbf{g},\mathbf{h}_{\mathbf{F}}}\bigr\|_{L^\infty(M)}.
\end{equation} In particular, if $\rho_M\leq1$, a section with bounded derivative is Hölder continuous with any exponent $0<\alpha<1$ for the pairs appearing in \eqref{eq:seminorma-holder-fibrada}.

Now let $\mathbf{F}_1,\mathbf{F}_2\to M$ be vector bundles with bundle metrics (Hermitian in the complex case) and compatible connections, $\mathbf{A}\in\Gamma^0(\mathbf{F}_1)$ and $\mathbf{B}\in\Gamma^0(\mathbf{F}_2)$. The product connection satisfies the Leibniz rule, so uniqueness of parallel transport implies \[
P^{\mathbf{F}_1\otimes \mathbf{F}_2}_{\gamma_{xy},1\to0}
=
P^{\mathbf{F}_1}_{\gamma_{xy},1\to0}
\otimes
P^{\mathbf{F}_2}_{\gamma_{xy},1\to0}.
\] Adding and subtracting $P^{\mathbf{F}_1}_{\gamma_{xy},1\to0}(\mathbf{A}(y))\otimes \mathbf{B}(x)$ gives \begin{equation}
\label{eq:producto-seminormas-holder-secciones}
[\mathbf{A}\otimes \mathbf{B}]_{\alpha,\mathbf{F}_1\otimes \mathbf{F}_2}
\leq
\|\mathbf{B}\|_{L^\infty(M,\mathbf{F}_2)}[\mathbf{A}]_{\alpha,\mathbf{F}_1}
+
\|\mathbf{A}\|_{L^\infty(M,\mathbf{F}_1)}[\mathbf{B}]_{\alpha,\mathbf{F}_2}.
\end{equation} The same estimate holds after any tensor contraction of uniformly bounded norm.

The iterated Leibniz rule distributes the $j$ derivatives between $\chi$ and $\mathbf{u}$. More precisely, for each $j\in\{0,\dots,k\}$ there are integers $N_{j,\ell}$, universal coefficients $C_{j,\ell,\nu}$, and permutations $\tau_{j,\ell,\nu}$ of the $j$ covariant indices such that \[
\nabla^j(\chi \mathbf{u})
=
\sum_{\ell=0}^{j}
\sum_{\nu=1}^{N_{j,\ell}}
C_{j,\ell,\nu}\,
\tau_{j,\ell,\nu}
\bigl(\nabla^{j-\ell}\chi\otimes\nabla^\ell \mathbf{u}\bigr).
\] Here a sum with $N_{j,\ell}=0$ is understood to be empty. Since $\chi$ has compact support, its covariant derivatives of every order, and the derivatives of those sections, are uniformly bounded. The preceding formula therefore controls the uniform norms of $\nabla^j(\chi \mathbf{u})$.

For the seminorm of $\nabla^k(\chi \mathbf{u})$, if $\ell<k$, the section $\nabla^\ell \mathbf{u}$ has bounded derivative $\nabla^{\ell+1}\mathbf{u}$. By \eqref{eq:seccion-c1-es-lipschitz-paralela-multiplicacion}, both $\nabla^\ell \mathbf{u}$ and $\nabla^{k-\ell}\chi$ are Hölder continuous with exponent $\alpha$. Estimate \eqref{eq:producto-seminormas-holder-secciones} controls the corresponding term. For $\ell=k$, apply the scalar version of the same estimate to the product $\chi\nabla^k\mathbf{u}$; the Hölder seminorm of $\nabla^k\mathbf{u}$ is part of \eqref{eq:norma-gamma-k-alpha}. Summing the terms gives the stated inequality. \end{proof}

We can now obtain the classical analogue of Lemma~\ref{lema:meyers-serrin-haz-E}. The proof is based on the triangular coordinate expression for covariant derivatives.

\begin{lemma}[Local characterization of Hölder sections] \label{lem:caracterizacion-local-holder-secciones} Let $(M,\mathbf{g})$ be a connected Riemannian manifold without boundary with $\operatorname{inj}(M)>0$. Let $(U,\phi)$ be a regular coordinate ball in $M$, and let $\mathbf{E}\longrightarrow M$ be a smooth vector bundle of rank $r$, equipped with a bundle metric $\mathbf{h}_{\mathbf{E}}$ and a compatible connection $\nabla^{\mathbf{E}}$. Let $(\mathbf{e}_1,\dots,\mathbf{e}_r)$ be a local frame of $\mathbf{E}$ defined on a neighborhood of $\overline U$, and write \[
\mathbf{u}=\sum_{a=1}^{r}u^a\mathbf{e}_a
\qquad\text{on $U$}.
\] Let $k\in\mathbb N_0$ and $0<\alpha<1$. Then \[
\mathbf{u}\in\Gamma^{k,\alpha}(\overline U,\mathbf{E})
\quad\Longleftrightarrow\quad
u^a\circ\phi^{-1}
\in
C^{k,\alpha}(\overline{\phi(U)})
\quad
\text{for every $a\in\{1,\dots,r\}$}.
\] Moreover, there are constants $c_{k,\alpha},C_{k,\alpha}>0$, independent of $\mathbf{u}$, such that \begin{equation}
\label{eq:equivalencia-local-holder-secciones}
\begin{aligned}
c_{k,\alpha}
\sum_{a=1}^{r}
\|u^a\circ\phi^{-1}\|_{C^{k,\alpha}(\overline{\phi(U)})}
&\leq
\|\mathbf{u}\|_{\Gamma^{k,\alpha}(\overline U,\mathbf{E})}\\
&\leq
C_{k,\alpha}
\sum_{a=1}^{r}
\|u^a\circ\phi^{-1}\|_{C^{k,\alpha}(\overline{\phi(U)})}.
\end{aligned}
\end{equation} \end{lemma}

\begin{proof} The case $k=0$ is Lemma~\ref{lem:comparacion-local-holder-transporte-componentes}. Suppose that $k\geq1$. Let $U'$ be a larger regular coordinate ball, so that $\overline U\subseteq U'$ and the chart and frame are defined on $U'$. In the chosen chart and frame, Lemma~\ref{lem:local-expression-higher-order} gives, for each $s\in\{1,\dots,k\}$, \begin{equation}
\label{eq:formula-triangular-holder-secciones}
(\nabla^s\mathbf{u})^a_{i_1\dots i_s}
=
\partial_{i_1}\cdots\partial_{i_s}u^a
+
\sum_{b=1}^{r}
\sum_{|\beta|\leq s-1}
(B_{s,\beta})^a_{b\,i_1\dots i_s}\partial^\beta u^b.
\end{equation} The coefficients $B_{s,\beta}$ are smooth on $U'$. Their derivatives of the orders appearing below are therefore uniformly bounded on each compact subset of $U'$.

We first prove the left inequality in \eqref{eq:equivalencia-local-holder-secciones}. Let $\widetilde{\mathbf{u}}\in\mathcal E_{\overline U}^{k,\alpha}(\mathbf{u})$. Choose a cutoff function $\chi\in C_c^\infty(U')$ equal to one on a neighborhood of $\overline U$ and set $\mathbf{v}:=\chi\widetilde{\mathbf{u}}$. If $\mathbf{v}=\displaystyle\sum_{a=1}^{r}v^a\mathbf{e}_a$ on $U'$, denote by $\widetilde v^a$ the extension of $v^a\circ\phi^{-1}$ by zero outside $\phi(U')$. Each $\widetilde v^a$ belongs to $C^k(\mathbb R^n)$ and extends $u^a\circ\phi^{-1}$.

Solving for the principal part in \eqref{eq:formula-triangular-holder-secciones} gives \begin{equation}
\label{eq:inversa-triangular-holder-secciones}
\partial_{i_1}\cdots\partial_{i_s}v^a
=
(\nabla^s\mathbf{v})^a_{i_1\dots i_s}
-
\sum_{b=1}^{r}
\sum_{|\beta|\leq s-1}
(B_{s,\beta})^a_{b\,i_1\dots i_s}\partial^\beta v^b.
\end{equation} Uniform equivalence of the fiber norms and induction on $s$ give \begin{equation}
\label{eq:cota-uniforme-componentes-desde-covariante-holder}
\sum_{a=1}^{r}\sum_{|\beta|\leq k}
\|D^\beta\widetilde v^a\|_{L^\infty(\mathbb R^n)}
\leq
C\sum_{s=0}^{k}\|\nabla^s\mathbf{v}\|_{L^\infty(M,\mathbf{E}_s)}.
\end{equation} Let us prove this estimate by induction. For $0\leq s\leq k$, set \[
S_s:=
\sum_{a=1}^{r}\sum_{|\beta|\leq s}
\|D^\beta\widetilde v^a\|_{L^\infty(\mathbb R^n)},
\qquad
T_s:=
\sum_{q=0}^{s}
\|\nabla^q\mathbf{v}\|_{L^\infty(M,\mathbf{E}_q)}.
\] Since $\mathbf{v}$ has compact support in the domain of the frame, extension by zero does not change the suprema of its components. Uniform equivalence of the norm of a section and the norms of its components gives $S_0\leq C_0T_0$.

Suppose that $S_{s-1}\leq C_{s-1}T_{s-1}$ for some $1\leq s\leq k$. Taking absolute values in \eqref{eq:inversa-triangular-holder-secciones}, using the uniform bounds on the coefficients $B_{s,\beta}$ on the compact set containing the support of $\mathbf{v}$, and summing over $a$ and all blocks $(i_1,\dots,i_s)$ gives \[
\sum_{a=1}^{r}\sum_{|\gamma|=s}
\|D^\gamma\widetilde v^a\|_{L^\infty(\mathbb R^n)}
\leq
C_s'\left(
\|\nabla^s\mathbf{v}\|_{L^\infty(M,\mathbf{E}_s)}+S_{s-1}
\right).
\] Here passing from ordered blocks of indices to multi-indices changes only the constant, since both families are finite. Consequently, \[
\begin{aligned}
S_s
&\leq
S_{s-1}+C_s'
\left(
\|\nabla^s\mathbf{v}\|_{L^\infty(M,\mathbf{E}_s)}+S_{s-1}
\right)\\
&\leq
\bigl((1+C_s')C_{s-1}+C_s'\bigr)T_s.
\end{aligned}
\] This proves the induction step. Taking the maximum of the constants appearing for $0\leq s\leq k$ yields \eqref{eq:cota-uniforme-componentes-desde-covariante-holder}.

It remains to control the Hölder seminorms of the derivatives of order $k$. The components of $\nabla^k\mathbf{v}$ are Hölder continuous in coordinates by Lemma~\ref{lem:comparacion-local-holder-transporte-componentes}, applied to the bundle $\mathbf{E}_k=T^{(0,k)}(TM)\otimes \mathbf{E}$ over an intermediate ball containing the support of $\chi$. If $|\beta|\leq k-1$, estimate \eqref{eq:cota-uniforme-componentes-desde-covariante-holder} controls the first derivatives of $D^\beta\widetilde v^b$. By the mean value theorem, these functions are Lipschitz and, since their supports lie in a fixed compact set, they are Hölder continuous with exponent $\alpha$ on $\mathbb R^n$. The inequality \[
[f\mathbf{g}]_{C^{0,\alpha}(K)}
\leq
\|f\|_{C^0(K)}[\mathbf{g}]_{C^{0,\alpha}(K)}
+
\|\mathbf{g}\|_{C^0(K)}[f]_{C^{0,\alpha}(K)}
\] applied to the lower-order terms in \eqref{eq:inversa-triangular-holder-secciones} gives \[
\sum_{a=1}^{r}
[\widetilde v^a]_{C^{k,\alpha}(\mathbb R^n)}
\leq
C\|\mathbf{v}\|_{\Gamma^{k,\alpha}(\mathbf{E})}.
\] Together with \eqref{eq:cota-uniforme-componentes-desde-covariante-holder}, this implies \[
\sum_{a=1}^{r}
\|u^a\circ\phi^{-1}\|_{C^{k,\alpha}(\overline{\phi(U)})}
\leq
C\|\widetilde{\mathbf{u}}\|_{\Gamma^{k,\alpha}(\mathbf{E})},
\] where we use Lemma~\ref{lem:multiplicacion-funcion-suave-holder-secciones} to control $\mathbf{v}=\chi\widetilde{\mathbf{u}}$. Taking the infimum over all extensions $\widetilde{\mathbf{u}}$ proves the first inequality in \eqref{eq:equivalencia-local-holder-secciones}.

We prove the reverse inequality. For each $a\in\{1,\dots,r\}$ choose an extension $w^a\in C^{k,\alpha}(\mathbb R^n)$ of $u^a\circ\phi^{-1}$. Let $\eta\in C_c^\infty(\phi(U'))$ be a cutoff function equal to one on a neighborhood of $\overline{\phi(U)}$. On $U'$ define \[
\mathbf{w}
:=
\sum_{a=1}^{r}\bigl((\eta w^a)\circ\phi\bigr)\mathbf{e}_a
\] and denote its extension by zero outside $U'$ by $\widetilde{\mathbf{w}}$. The presence of $\eta$ ensures that $\widetilde{\mathbf{w}}$ is a global section of class $C^k$ and that $\widetilde{\mathbf{w}}\restriction_U=\mathbf{u}$.

The direct formula \eqref{eq:formula-triangular-holder-secciones}, the Leibniz rule, and Lemma~\ref{lem:comparacion-local-holder-transporte-componentes}, again applied to the bundle $\mathbf{E}_k$, give \[
\|\widetilde{\mathbf{w}}\|_{\Gamma^{k,\alpha}(\mathbf{E})}
\leq
C\sum_{a=1}^{r}\|w^a\|_{C^{k,\alpha}(\mathbb R^n)}.
\] The terms of order $k$ are controlled by the seminorms of the $w^a$. The lower-order terms have an additional bounded derivative and are controlled by the mean value theorem, exactly as in the first part. Taking the infimum first over the scalar extensions $w^a$ and then in the definition of the norm of $\mathbf{u}$ gives the second inequality in \eqref{eq:equivalencia-local-holder-secciones}. The two constructions also prove equivalence of membership. \end{proof}

The following localized version of the norm will allow us to combine estimates obtained in different charts.

\begin{proposition}[Equivalence of the localized Hölder norm] \label{prop:equivalencia-norma-holder-localizada-haces} Let $(M,\mathbf{g})$ be a compact connected Riemannian manifold without boundary. Let $\mathbf{E}\longrightarrow M$ be a smooth vector bundle equipped with a bundle metric and a compatible connection. Consider regular coordinate balls $(U_i,\phi_i)_{i=1}^{N}$ covering $M$ and smooth cutoff functions $(\psi_i)_{i=1}^{N}$ such that \[
0\leq\psi_i\leq1,
\qquad
\operatorname{supp}(\psi_i)\Subset U_i,
\qquad
\sum_{i=1}^{N}\psi_i=1.
\] For each $i$, let $(\mathbf{e}_{i,1},\dots,\mathbf{e}_{i,r})$ be a local frame of $\mathbf{E}$ defined on a neighborhood of $\overline{U_i}$, and write \[
\psi_i \mathbf{u}
=
\sum_{a=1}^{r}u_i^a\mathbf{e}_{i,a}
\qquad\text{on $U_i$}.
\] If $k\in\mathbb N_0$ and $0<\alpha<1$, there are constants $c,C>0$ such that \begin{equation}
\label{eq:norma-holder-localizada-haces}
c\|\mathbf{u}\|_{\Gamma^{k,\alpha}(\mathbf{E})}
\leq
\sum_{i=1}^{N}\sum_{a=1}^{r}
\|u_i^a\circ\phi_i^{-1}\|_{C^{k,\alpha}(\overline{\phi_i(U_i)})}
\leq
C\|\mathbf{u}\|_{\Gamma^{k,\alpha}(\mathbf{E})}
\end{equation} for every $\mathbf{u}\in\Gamma^{k,\alpha}(\mathbf{E})$. \end{proposition}

\begin{proof} Lemma~\ref{lem:multiplicacion-funcion-suave-holder-secciones} gives a constant $C_0>0$, independent of $i$, such that \[
\|\psi_i \mathbf{u}\|_{\Gamma^{k,\alpha}(\mathbf{E})}
\leq
C_0\|\mathbf{u}\|_{\Gamma^{k,\alpha}(\mathbf{E})}.
\] Uniformity follows because only a finite family of cutoff functions is involved. Lemma~\ref{lem:caracterizacion-local-holder-secciones}, applied to $\psi_i \mathbf{u}\restriction_{U_i}$, gives the second inequality in \eqref{eq:norma-holder-localizada-haces}.

We prove the first. Set $K_i:=\phi_i(\operatorname{supp}\psi_i)\Subset\phi_i(U_i)$ and choose $\eta_i\in C_c^\infty(\phi_i(U_i))$ such that $\eta_i=1$ on a neighborhood of $K_i$. For each $a$ take an extension $v_i^a\in C^{k,\alpha}(\mathbb R^n)$ of $u_i^a\circ\phi_i^{-1}$. The function $\eta_i v_i^a$, extended by zero outside $\phi_i(U_i)$, remains an extension of that component. On $U_i$ reconstruct the section \[
\sum_{a=1}^{r}\bigl((\eta_i v_i^a)\circ\phi_i\bigr)\mathbf{e}_{i,a}
\] and denote its extension by zero outside $U_i$ by $\widetilde{\mathbf{u}}_i$. By construction, $\widetilde{\mathbf{u}}_i=\psi_i \mathbf{u}$. The local comparison and product rule give \[
\|\psi_i \mathbf{u}\|_{\Gamma^{k,\alpha}(\mathbf{E})}
\leq
C_i\sum_{a=1}^{r}\|v_i^a\|_{C^{k,\alpha}(\mathbb R^n)}.
\] Taking the infimum over the extensions $v_i^a$, summing over $i$, and using $\mathbf{u}=\displaystyle\sum_{i=1}^{N}\psi_i \mathbf{u}$ gives \[
\|\mathbf{u}\|_{\Gamma^{k,\alpha}(\mathbf{E})}
\leq
C\sum_{i=1}^{N}\sum_{a=1}^{r}
\|u_i^a\circ\phi_i^{-1}\|_{C^{k,\alpha}(\overline{\phi_i(U_i)})}.
\] This completes the proof. \end{proof}

\begin{corollary}[Completeness of spaces of Hölder sections] \label{cor:completitud-gamma-k-alpha} Let $(M,\mathbf{g})$ be a compact connected Riemannian manifold without boundary, and let $\mathbf{E}\to M$ be a smooth vector bundle with a bundle metric and a compatible connection. Let $k\in\mathbb N_0$ and $0<\alpha<1$. Then $\Gamma^{k,\alpha}(\mathbf{E})$ is a Banach space. \end{corollary}

\begin{proof} Let $(\mathbf{u}_j)_{j\in\mathbb N}$ be a Cauchy sequence in $\Gamma^{k,\alpha}(\mathbf{E})$. For each $i\in\{1,\dots,N\}$ write \[
\psi_i \mathbf{u}_j
=
\sum_{a=1}^{r}u_{i,j}^a \mathbf{e}_{i,a}
\qquad\text{on $U_i$}.
\] For each pair $(i,a)$, the components of $\psi_i \mathbf{u}_j$ form a Cauchy sequence in $C^{k,\alpha}(\overline{\phi_i(U_i)})$ by \eqref{eq:norma-holder-localizada-haces}. Proposition~\ref{prop:completitud-holder-euclidiano} gives a limit $f_i^a$ in that space. Since the extension norm controls the uniform norm on the closure, $f_i^a$ vanishes outside the fixed compact set $K_i=\phi_i(\operatorname{supp}\psi_i)$, as do all the components of the sequence.

Choose global extensions of the $f_i^a$, multiply them by the functions $\eta_i$ from the preceding proof, and reconstruct global sections $\mathbf{v}_i$ supported in $U_i$. The estimate obtained in that proof gives \[
\|\psi_i \mathbf{u}_j-\mathbf{v}_i\|_{\Gamma^{k,\alpha}(\mathbf{E})}
\leq
C_i\sum_{a=1}^{r}
\|(u_{i,j}^a\circ\phi_i^{-1})-f_i^a\|_{C^{k,\alpha}(\overline{\phi_i(U_i)})},
\] and the right-hand side converges to zero. Define \[
\mathbf{u}:=\sum_{i=1}^{N}\mathbf{v}_i.
\] Since $\mathbf{u}_j=\displaystyle\sum_{i=1}^{N}\psi_i \mathbf{u}_j$ and the sum is finite, \[
\|\mathbf{u}_j-\mathbf{u}\|_{\Gamma^{k,\alpha}(\mathbf{E})}
\leq
\sum_{i=1}^{N}
\|\psi_i \mathbf{u}_j-\mathbf{v}_i\|_{\Gamma^{k,\alpha}(\mathbf{E})}
\longrightarrow0.
\] Therefore, $\mathbf{u}\in\Gamma^{k,\alpha}(\mathbf{E})$ and the space is complete. \end{proof}

We can now transfer the Euclidean embeddings to bundle sections.

\begin{theorem}[Sobolev--Morrey embedding for vector bundles] \label{teo:sobolev-morrey-haces} \index{Sobolev Morrey theorem@Sobolev--Morrey theorem!for vector bundles} Let $(M,\mathbf{g})$ be a compact Riemannian manifold without boundary of dimension $n$, and let $\mathbf{E}\longrightarrow M$ be a smooth vector bundle equipped with a bundle metric $\mathbf{h}_{\mathbf{E}}$ and a compatible connection $\nabla^{\mathbf{E}}$. Let $m\in\mathbb N$, $1\leq p<\infty$, $k\in\mathbb N_0$, and $0<\alpha<1$. Suppose that \begin{equation}
\label{eq:condicion-sobolev-morrey-haces}
k+\alpha
\leq
m-\frac np.
\end{equation} Then each element of $W^{m,p}(M,\mathbf{E})$ admits a unique representative in $\Gamma^{k,\alpha}(\mathbf{E})$, and the map assigning this representative defines a continuous linear embedding \[
W^{m,p}(M,\mathbf{E})
\hookrightarrow
\Gamma^{k,\alpha}(\mathbf{E}).
\] In particular, there is a constant $C>0$ such that \begin{equation}
\label{eq:estimacion-sobolev-morrey-haces}
\|\mathbf{u}\|_{\Gamma^{k,\alpha}(\mathbf{E})}
\leq
C\|\mathbf{u}\|_{W^{m,p}(M,\mathbf{E})}
\end{equation} for every $\mathbf{u}\in W^{m,p}(M,\mathbf{E})$. \end{theorem}

\begin{proof} We may work on each connected component of $M$, and therefore assume that $M$ is connected. Choose regular coordinate balls, cutoff functions, and local frames as in Proposition~\ref{prop:equivalencia-norma-holder-localizada-haces}. For $\mathbf{u}\in W^{m,p}(M,\mathbf{E})$ write \[
\psi_i \mathbf{u}
=
\sum_{a=1}^{r}u_i^a\mathbf{e}_{i,a}
\qquad\text{on $U_i$}.
\] The Leibniz estimate used in the proof of Theorem~\ref{teo: encaje de sobolev haces}, in particular estimate \eqref{eq: acotacion particion de la unidad y norma sobolev auxiliar para encaje de Sobolev y Rellich--Kondrashov} with $j=0$, gives a constant $C_0>0$ such that \begin{equation}
\label{eq:multiplicacion-cortes-sobolev-morrey-haces}
\|\psi_i \mathbf{u}\|_{W^{m,p}(M,\mathbf{E})}
\leq
C_0\|\mathbf{u}\|_{W^{m,p}(M,\mathbf{E})}
\end{equation} for every $i\in\{1,\dots,N\}$.

Lemma~\ref{lema:meyers-serrin-haz-E} shows that $u_i^a\circ\phi_i^{-1}\in W^{m,p}(\phi_i(U_i))$ and that \begin{equation}
\label{eq:cota-componentes-sobolev-morrey-haces}
\sum_{a=1}^{r}
\|u_i^a\circ\phi_i^{-1}\|_{W^{m,p}(\phi_i(U_i))}
\leq
C_i
\|\psi_i \mathbf{u}\|_{W^{m,p}(U_i,\mathbf{E}\restriction_{U_i})}.
\end{equation} Since $\operatorname{supp}(\psi_i)\Subset U_i$, each component has compact support contained in $\phi_i(U_i)$ and therefore belongs to $W_0^{m,p}(\phi_i(U_i))$. Theorem~\ref{teo:sobolev-morrey-euclidiano} and \eqref{eq:condicion-sobolev-morrey-haces} give \[
\|u_i^a\circ\phi_i^{-1}\|_{C^{k,\alpha}(\overline{\phi_i(U_i)})}
\leq
C_i'
\|u_i^a\circ\phi_i^{-1}\|_{W^{m,p}(\phi_i(U_i))}.
\]

Take the Euclidean representatives of the components, reconstruct the section $\psi_i \mathbf{u}$ on $U_i$, and extend it by zero outside $U_i$. The sum of these sections represents $\mathbf{u}$, since it agrees almost everywhere with $\displaystyle\sum_{i=1}^{N}\psi_i \mathbf{u}=\mathbf{u}$. Proposition~\ref{prop:equivalencia-norma-holder-localizada-haces}, the preceding estimates, and finiteness of the cover give \[
\begin{aligned}
\|\mathbf{u}\|_{\Gamma^{k,\alpha}(\mathbf{E})}
&\leq
C
\sum_{i=1}^{N}
\sum_{a=1}^{r}
\|u_i^a\circ\phi_i^{-1}\|_{C^{k,\alpha}(\overline{\phi_i(U_i)})}\\
&\leq
C
\sum_{i=1}^{N}
\sum_{a=1}^{r}
\|u_i^a\circ\phi_i^{-1}\|_{W^{m,p}(\phi_i(U_i))}\\
&\leq
C
\sum_{i=1}^{N}
\|\psi_i \mathbf{u}\|_{W^{m,p}(M,\mathbf{E})}
\leq
C'\|\mathbf{u}\|_{W^{m,p}(M,\mathbf{E})}.
\end{aligned}
\] This proves \eqref{eq:estimacion-sobolev-morrey-haces}.

Uniqueness of the representative follows as in the Euclidean case. If a continuous section vanishes almost everywhere but is nonzero at some point, its norm is positive on an open neighborhood of that point, which has positive Riemann--Lebesgue measure. Thus the section must vanish throughout $M$. \end{proof}

\begin{theorem}[Rellich--Kondrashov into Hölder sections] \label{teo:rellich-kondrashov-holder-haces} \index{Rellich Kondrashov theorem@Rellich--Kondrashov theorem!into Hölder sections} Let $(M,\mathbf{g})$ be a compact Riemannian manifold without boundary of dimension $n$, and let $\mathbf{E}\to M$ be a smooth vector bundle with a bundle metric and a compatible connection. Let $m\in\mathbb N$, $1\leq p<\infty$, $k\in\mathbb N_0$, and $0<\alpha<1$. If \begin{equation}
\label{eq:condicion-compacta-holder-haces}
k+\alpha
<
m-\frac np,
\end{equation} then the embedding \[
W^{m,p}(M,\mathbf{E})
\hookrightarrow
\Gamma^{k,\alpha}(\mathbf{E})
\] is compact. \end{theorem}

\begin{proof} Let $(\mathbf{u}_j)_{j\in\mathbb N}$ be a bounded sequence in $W^{m,p}(M,\mathbf{E})$. Use the cover, cutoff functions, and frames from the preceding proof. For each $i\in\{1,\dots,N\}$ write \[
\psi_i \mathbf{u}_j
=
\sum_{a=1}^{r}u_{i,j}^a\mathbf{e}_{i,a}
\qquad\text{on $U_i$}.
\] Estimates \eqref{eq:multiplicacion-cortes-sobolev-morrey-haces} and \eqref{eq:cota-componentes-sobolev-morrey-haces} show that, for each pair $(i,a)$, the sequence \[
\bigl(u_{i,j}^a\circ\phi_i^{-1}\bigr)_{j\in\mathbb N}
\] is bounded in $W_0^{m,p}(\phi_i(U_i))$. By Theorem~\ref{teo:rellich-kondrashov-holder-euclidiano}, each of these sequences has a subsequence converging in $C^{k,\alpha}(\overline{\phi_i(U_i)})$.

Only the $Nr$ pairs of indices are involved. Extracting successive subsequences and relabeling the last one, we may assume that \[
u_{i,j}^a\circ\phi_i^{-1}
\quad\text{converges in}
\quad
C^{k,\alpha}(\overline{\phi_i(U_i)})
\] for every $i\in\{1,\dots,N\}$ and every $a\in\{1,\dots,r\}$. In particular, the localized components form Cauchy sequences in those spaces. The first inequality in \eqref{eq:norma-holder-localizada-haces}, applied to $\mathbf{u}_j-\mathbf{u}_\ell$, gives \[
\|\mathbf{u}_j-\mathbf{u}_\ell\|_{\Gamma^{k,\alpha}(\mathbf{E})}
\leq
C
\sum_{i=1}^{N}
\sum_{a=1}^{r}
\|(u_{i,j}^a-u_{i,\ell}^a)\circ\phi_i^{-1}\|_{C^{k,\alpha}(\overline{\phi_i(U_i)})}.
\] The right-hand side converges to zero as $j,\ell\to\infty$. Thus the subsequence is Cauchy in $\Gamma^{k,\alpha}(\mathbf{E})$ and converges in that space by Corollary~\ref{cor:completitud-gamma-k-alpha}. This proves compactness of the embedding. \end{proof}

\begin{corollary}[Regularity exponents for sections] \label{cor:exponentes-regulares-holder-haces} Let $(M,\mathbf{g})$ be a compact Riemannian manifold without boundary of dimension $n$, and let $\mathbf{E}\to M$ be a smooth vector bundle with a bundle metric and a compatible connection. Let $m\in\mathbb N$ and $1\leq p<\infty$. Let \[
\sigma:=m-\frac np>0.
\] Then the following assertions hold. \begin{enumerate}[label=(\alph*)] \item If $\sigma=\ell+\theta$, with $\ell\in\mathbb N_0$ and $0<\theta<1$, then \[
W^{m,p}(M,\mathbf{E})
\hookrightarrow
\Gamma^{\ell,\theta}(\mathbf{E})
\] continuously, and the embedding into $\Gamma^{\ell,\beta}(\mathbf{E})$ is compact for every $0<\beta<\theta$. \item If $\sigma=s\in\mathbb N$, then for every $0<\alpha<1$ the embedding \[
W^{m,p}(M,\mathbf{E})
\hookrightarrow
\Gamma^{s-1,\alpha}(\mathbf{E})
\] is compact and, in particular, continuous. \end{enumerate} \end{corollary} \begin{proof} If \(\sigma=\ell+\theta\) with \(0<\theta<1\), Theorem~\ref{teo:sobolev-morrey-haces} applied to \(k=\ell\) and \(\alpha=\theta\) gives the continuous embedding, since \[
\ell+\theta=m-\frac np.
\] If \(0<\beta<\theta\), then \[
\ell+\beta<m-\frac np,
\] and Theorem~\ref{teo:rellich-kondrashov-holder-haces} gives compactness of the embedding into \(\Gamma^{\ell,\beta}(\mathbf E)\).

If \(\sigma=s\in\mathbb N\) and \(0<\alpha<1\), we have \[
(s-1)+\alpha<s=m-\frac np.
\] The same Rellich--Kondrashov theorem, with \(k=s-1\), gives the compact embedding into \(\Gamma^{s-1,\alpha}(\mathbf E)\). Every compact linear operator between normed spaces is continuous, which also proves the final assertion. \end{proof}

\subsection{Independence of Sobolev spaces in the compact case} As in the setting of Riemannian manifolds, Sobolev spaces on vector bundles over a compact Riemannian manifold are independent of the choice of metric and connection. To prove this, we establish some useful preliminary results. In particular, we formulate an analogue of Lemma~\ref{lema dificil de probar} in the vector bundle setting. \begin{lemma}\label{lema: independencia de normas en compactos} Let $M$ be a smooth manifold with or without boundary, let $\mathbf{g}$ and $\widetilde{\mathbf{g}}$ be two Riemannian metrics on $M$, let $\nabla^{TM}$ and $\widetilde{\nabla}^{TM}$ be two connections on $TM$, let $\mathbf{E}\longrightarrow M$ be a smooth vector bundle of finite rank, let $\mathbf{h}$ and $\widetilde{\mathbf{h}}$ be two smooth bundle metrics on $\mathbf{E}$, and let $\nabla^{\mathbf{E}}$ and $\widetilde{\nabla}^{\mathbf{E}}$ be two connections on $\mathbf{E}$.

Let $K\subseteq M$ be compact and $m\geq 0$. Then there are constants $c_{m},C_{m}>0$ such that \[
c_{m}\displaystyle\sum_{j=0}^{m}|\nabla^{j}\mathbf{u}(x)|^{2}_{\mathbf{g},\mathbf{h}}
\leq
\displaystyle\sum_{j=0}^{m}|\widetilde{\nabla}^{j}\mathbf{u}(x)|^{2}_{\widetilde{\mathbf{g}},\widetilde{\mathbf{h}}}
\leq
C_{m}\displaystyle\sum_{j=0}^{m}|\nabla^{j}\mathbf{u}(x)|^{2}_{\mathbf{g},\mathbf{h}},
\qquad \forall x\in K,
\] for every section $\mathbf{u}\in \Gamma^{m}(\mathbf{E})$. \end{lemma} \begin{proof} If $K=\varnothing$, the assertion is immediate. Suppose that $K\neq\varnothing$. The connections $\nabla^{TM}$ and $\nabla^{\mathbf E}$ induce connections on $T^{(0,j)}(TM)\otimes\mathbf E$, also denoted by $\nabla$. By Definition~\ref{suma de whitney}, we may consider \[
 \mathscr T_m(\mathbf E)
 :=\bigoplus_{j=0}^{m}T^{(0,j)}(TM)\otimes\mathbf E.
\] For $\mathbf u\in\Gamma(\mathbf E)$ write \[
 J_\nabla^m\mathbf u
 :=(\mathbf u,\nabla\mathbf u,\dots,\nabla^m\mathbf u)
\] and, for each $p\in M$, define \[
 \mathscr D_{\nabla,p}(\mathbf E)
 :=\{J_\nabla^m\mathbf u(p)\mid\mathbf u\in\Gamma(\mathbf E)\}
 \subseteq\mathscr T_m(\mathbf E)_p.
\] Each $\mathscr D_{\nabla,p}(\mathbf E)$ is a vector subspace by linearity of covariant derivatives.

Let $(U,\phi)$ be a smooth chart, let $(\mathbf e_1,\dots,\mathbf e_r)$ be a smooth local frame of $\mathbf E$ on $U$, let $p\in U$, and set $x:=\phi(p)$ and $N:={n+m\choose m}$. If $\mathbf u=\displaystyle\sum_{a=1}^{r}u^a\mathbf e_a$, define the linear map \[
 F_{\nabla,p}^{\phi,\mathbf e}\colon
 \mathscr D_{\nabla,p}(\mathbf E)\longrightarrow\mathbb K^{rN},
 \qquad
 F_{\nabla,p}^{\phi,\mathbf e}(J_\nabla^m\mathbf u(p))
 :=\bigl(D^\alpha(u^a\circ\phi^{-1})(x)\bigr)_{
 1\leq a\leq r,\ |\alpha|\leq m},
\] using the multi-index notation from Definition~\ref{def:nociones-fundamentales-multiindice}. Lemma~\ref{lem:local-expression-higher-order}, in particular formula~\eqref{eq:higher-order-formula}, gives \begin{equation}
 (\nabla^j\mathbf u)^a_{i_1\dots i_j}
 =\partial_{i_1}\cdots\partial_{i_j}(u^a\circ\phi^{-1})
 +\sum_{b=1}^{r}\sum_{|\alpha|<j}
 A^{ab,\alpha}_{i_1\dots i_j}
 D^\alpha(u^b\circ\phi^{-1}),
 \label{eq:jet-covariante-triangular-secciones}
\end{equation} with smooth coefficients on $U$. Induction on $j$ in \eqref{eq:jet-covariante-triangular-secciones} shows that the covariant jet determines all partial derivatives of the components through order $m$. Consequently, $F_{\nabla,p}^{\phi,\mathbf e}$ is well defined and injective.

Let $a=(a_\alpha^b)_{1\leq b\leq r,\ |\alpha|\leq m}\in\mathbb K^{rN}$. Proposition~\ref{funciones flan} gives a function $\chi\in C_c^\infty(U)$ equal to one on a neighborhood of $p$. The section \[
 \mathbf v(q):=
 \chi(q)\sum_{b=1}^{r}
 \left(\sum_{|\alpha|\leq m}
 \frac{a_\alpha^b}{\alpha!}(\phi(q)-x)^\alpha\right)\mathbf e_b(q),
 \qquad q\in U,
\] extended by zero outside $U$ belongs to $\Gamma(\mathbf E)$ and satisfies $D^\alpha(v^b\circ\phi^{-1})(x)=a_\alpha^b$. Thus $F_{\nabla,p}^{\phi,\mathbf e}$ is surjective and hence an isomorphism.

Let $(e_1,\dots,e_{rN})$ be the standard basis of $\mathbb K^{rN}$. Expressing the covariant jets in terms of the prescribed partial jets using \eqref{eq:jet-covariante-triangular-secciones} shows that the components of the local sections \[
 \boldsymbol\sigma_\nu(q)
 :=(F_{\nabla,q}^{\phi,\mathbf e})^{-1}(e_\nu),
 \qquad q\in U,\quad \nu\in\{1,\dots,rN\},
\] are smooth. Their values form a basis of $\mathscr D_{\nabla,q}(\mathbf E)$ for every $q\in U$. By Lemma~\ref{criterio del marco local suave para subhaces}, \[
 \mathscr D_\nabla(\mathbf E)
 :=\bigcup_{p\in M}\mathscr D_{\nabla,p}(\mathbf E)
\] is a smooth vector subbundle of $\mathscr T_m(\mathbf E)$ of rank $rN$. The connections $\widetilde\nabla^{TM}$ and $\widetilde\nabla^{\mathbf E}$ similarly produce the subbundle $\mathscr D_{\widetilde\nabla}(\mathbf E)$ and the isomorphisms $F_{\widetilde\nabla,p}^{\phi,\mathbf e}$.

The compact set $K$ is covered by finitely many compact regular coordinate balls and, if $M$ has boundary, compact regular coordinate half-balls \[
 \overline B_1,\dots,\overline B_\ell,
\] each contained in the domain of a chart $\phi_i$ and a local frame $\mathbf e_i=(\mathbf e_{1,i},\dots,\mathbf e_{r,i})$ of $\mathbf E$. For $p\in\overline B_i$ and $v\in\mathbb K^{rN}$ set \[
 \begin{aligned}
 Q_i(p,v)
 &:=\displaystyle\sum_{j=0}^{m}
 \left|\bigl((F_{\nabla,p}^{\phi_i,\mathbf e_i})^{-1}(v)\bigr)_j
 \right|_{\mathbf g,\mathbf h}^{2},\\
 \widetilde Q_i(p,v)
 &:=\displaystyle\sum_{j=0}^{m}
 \left|\bigl((F_{\widetilde\nabla,p}^{\phi_i,\mathbf e_i})^{-1}(v)\bigr)_j
 \right|_{\widetilde{\mathbf g},\widetilde{\mathbf h}}^{2}.
 \end{aligned}
\] The functions $Q_i$ and $\widetilde Q_i$ are continuous and, for each $p$, are positive definite quadratic forms in $v$. If \[
 \mathbb S_{\mathbb K}^{rN}
 :=\{v\in\mathbb K^{rN}\mid\|v\|_2=1\},
\] then \[
 R_i(p,v):=\frac{\widetilde Q_i(p,v)}{Q_i(p,v)}
\] is continuous and strictly positive on the compact set $\overline B_i\times\mathbb S_{\mathbb K}^{rN}$. Define \[
 c_{i,m}:=\min_{(p,v)\in\overline B_i\times\mathbb S_{\mathbb K}^{rN}}
 R_i(p,v)>0,
 \qquad
 C_{i,m}:=\max_{(p,v)\in\overline B_i\times\mathbb S_{\mathbb K}^{rN}}
 R_i(p,v)<\infty.
\] Homogeneity of $Q_i$ and $\widetilde Q_i$ implies that \begin{equation}
 c_{i,m}Q_i(p,v)
 \leq\widetilde Q_i(p,v)
 \leq C_{i,m}Q_i(p,v)
 \label{eq:comparacion-local-formas-cuadraticas-jets-secciones}
\end{equation} for every $p\in\overline B_i$ and every $v\in\mathbb K^{rN}$.

Let $\mathbf u\in\Gamma^m(\mathbf E)$ and let $p\in\overline B_i$. The preceding polynomial construction gives a smooth section whose partial derivative jet at $p$ agrees with that of $\mathbf u$. Formula~\eqref{eq:higher-order-formula} then shows that the covariant jets of the two sections agree at $p$ for either pair of connections. Thus \[
 \begin{aligned}
 F_{\nabla,p}^{\phi_i,\mathbf e_i}(J_\nabla^m\mathbf u(p))
 &=\bigl(D^\alpha(u_i^a\circ\phi_i^{-1})(\phi_i(p))\bigr)_{
 1\leq a\leq r,\ |\alpha|\leq m}\\
 &=F_{\widetilde\nabla,p}^{\phi_i,\mathbf e_i}
 (J_{\widetilde\nabla}^m\mathbf u(p)).
 \end{aligned}
\] Applying \eqref{eq:comparacion-local-formas-cuadraticas-jets-secciones} to this vector gives \[
 c_{i,m}\sum_{j=0}^{m}|\nabla^j\mathbf u(p)|_{\mathbf g,\mathbf h}^{2}
 \leq
 \sum_{j=0}^{m}|\widetilde\nabla^j\mathbf u(p)|_{
 \widetilde{\mathbf g},\widetilde{\mathbf h}}^{2}
 \leq
 C_{i,m}\sum_{j=0}^{m}|\nabla^j\mathbf u(p)|_{\mathbf g,\mathbf h}^{2}.
\] Taking \[
 c_m:=\min_{1\leq i\leq\ell}c_{i,m}>0,
 \qquad
 C_m:=\max_{1\leq i\leq\ell}C_{i,m}<\infty,
\] yields the inequality throughout $K$. \end{proof}

\section{Interpolation: the Gagliardo--Nirenberg inequality} In this section we study the \textit{Gagliardo--Nirenberg interpolation} inequality on manifolds (and in vector bundles). The inequality states the following: \begin{theorem}[Gagliardo--Nirenberg inequality for vector bundles]\label{teo:gagliardo-nirenberg-haces}\index{Gagliardo Nirenberg inequality for vector bundles@Gagliardo--Nirenberg inequality for vector bundles} Let \((M,\mathbf{g})\) be a compact Riemannian manifold without boundary of dimension \(n\), and let \(\mathbf{E}\longrightarrow M\) be a smooth vector bundle of finite rank, equipped with a bundle metric \(\mathbf{h}_{\mathbf{E}}\) and a compatible connection \(\nabla^{\mathbf{E}}\). Let \(m\geq 1\), \(0\leq j\leq m\), and \(1\leq p<\infty\). Then there is a constant \(C>0\), depending only on \(M,\mathbf{g},\mathbf{E},\mathbf{h}_{\mathbf{E}},\nabla^{\mathbf{E}},m,j\) and \(p\), such that \[
\|\nabla^{j}\mathbf{F}\|_{L^p(M,T^{(0,j)}(TM)\otimes \mathbf{E})}
\leq
C\,
\|\mathbf{F}\|_{L^p(M,\mathbf{E})}^{1-\frac{j}{m}}
\left(
\|\mathbf{F}\|_{L^p(M,\mathbf{E})}
+
\|\nabla^{m}\mathbf{F}\|_{L^p(M,T^{(0,m)}(TM)\otimes \mathbf{E})}
\right)^{\frac{j}{m}},
\] for every \(\mathbf{F}\in W^{m,p}(M,\mathbf{E})\). Here we interpret \(\nabla^0\mathbf{F}=\mathbf{F}\), so that \(T^{(0,0)}(TM)\otimes \mathbf{E}=\mathbf{E}\). \end{theorem}

\begin{proof} The cases \(j=0\) and \(j=m\) are immediate. If \(j=0\), the inequality reduces to \(\|\mathbf{F}\|_{L^p(M,\mathbf{E})}\leq C\|\mathbf{F}\|_{L^p(M,\mathbf{E})}\); if \(j=m\), it reduces to \[
\|\nabla^m\mathbf{F}\|_{L^p(M,T^{(0,m)}(TM)\otimes \mathbf{E})}
\leq
C
\left(
\|\mathbf{F}\|_{L^p(M,\mathbf{E})}
+
\|\nabla^m\mathbf{F}\|_{L^p(M,T^{(0,m)}(TM)\otimes \mathbf{E})}
\right).
\] Thus suppose that \(0<j<m\).

We first prove the inequality for smooth sections. Let \(\mathbf{F}\in\Gamma(\mathbf{E})\). Since \(M\) is compact, we may take a finite cover \((U_{\alpha},\phi_{\alpha})_{\alpha=1}^{N}\) of \(M\) by regular coordinate balls such that $\overline{U_\alpha}$ lies in the domain of a smooth local frame \((\mathbf{e}_{1,\alpha},\dots,\mathbf{e}_{r,\alpha})\) of \(\mathbf{E}\). Also take a partition of unity \((\psi_{\alpha})_{\alpha=1}^{N}\) subordinate to this cover, with \(\supp(\psi_{\alpha})\subseteq U_{\alpha}\) for every \(\alpha\).

For each \(\alpha\in\{1,\dots,N\}\) write \[
(\psi_{\alpha}\mathbf{F})\restriction_{U_{\alpha}}
=
\displaystyle\sum_{a=1}^{r}F_{\alpha}^{a}\mathbf{e}_{a,\alpha},
\] and denote \(\Omega_{\alpha}:=\phi_{\alpha}(U_{\alpha})\). Since \(\supp(\psi_{\alpha})\subseteq U_{\alpha}\), we have \(\supp(F_{\alpha}^{a}\circ\phi_{\alpha}^{-1})\Subset\Omega_{\alpha}\) for every \(a\in\{1,\dots,r\}\).

Define the extension of \(F_{\alpha}^{a}\circ\phi_{\alpha}^{-1}\) by zero to \(\mathbb R^n\) by \[
\widetilde{F_{\alpha}^{a}\circ\phi_{\alpha}^{-1}}(y)
:=
\begin{cases}
F_{\alpha}^{a}\circ\phi_{\alpha}^{-1}(y), & y\in\Omega_{\alpha},\\
0, & y\in\mathbb R^n\setminus\Omega_{\alpha}.
\end{cases}
\] Since \(F_{\alpha}^{a}\circ\phi_{\alpha}^{-1}\) is smooth and has compact support contained in \(\Omega_{\alpha}\), its extension by zero belongs to \(C_c^\infty(\mathbb R^n)\). Moreover, for every multi-index \(\beta\), we have \[
D^\beta\widetilde{F_{\alpha}^{a}\circ\phi_{\alpha}^{-1}}
=
\widetilde{D^\beta(F_{\alpha}^{a}\circ\phi_{\alpha}^{-1})}.
\] Consequently, for each integer \(s\geq 0\), \[
\displaystyle\sum_{|\beta|=s}
\|D^\beta\widetilde{F_{\alpha}^{a}\circ\phi_{\alpha}^{-1}}\|_{L^p(\mathbb R^n)}
=
\displaystyle\sum_{|\beta|=s}
\|D^\beta(F_{\alpha}^{a}\circ\phi_{\alpha}^{-1})\|_{L^p(\Omega_{\alpha})}.
\]

We first establish the Euclidean inequality used here. For
$v\in C_c^\infty(\mathbb R^n)$, Taylor's formula along a coordinate line
and translation invariance of the $L^p(\mathbb R^n)$ norm give, for every $h>0$,
\[
\|\partial_i v\|_{L^p(\mathbb R^n)}
\leq\frac2h\|v\|_{L^p(\mathbb R^n)}+\frac h2\|\partial_i^2v\|_{L^p(\mathbb R^n)}.
\]
Minimizing in $h$ gives
$\|\partial_i v\|_{L^p(\mathbb R^n)}\leq2\|v\|_{L^p(\mathbb R^n)}^{1/2}
\|\partial_i^2v\|_{L^p(\mathbb R^n)}^{1/2}$. Put
$A_k=\sum_{|\alpha|=k}\|D^\alpha v\|_{L^p(\mathbb R^n)}$ for $0\leq k\leq m$.
Apply this estimate to $D^{\alpha-e_i}v$, choosing an index with
$\alpha_i>0$, and sum over $|\alpha|=k$. The finite Cauchy--Schwarz
inequality yields
\[
A_k\leq C_{n,m}(A_{k-1}A_{k+1})^{1/2},\qquad 0<k<m.
\]
Take $C_{n,m}\geq1$. The sequence
$\log A_k-k(m-k)\log C_{n,m}$ is convex, so comparison with the chord
joining its endpoint values proves
$A_j\leq C_{n,m}^{j(m-j)}A_0^{1-j/m}A_m^{j/m}$.
If some $A_k$ vanishes, apply the same argument to $A_k+\tau$ and
let $\tau\to0^{+}$; the preceding inequality remains valid after
this replacement. This proves the required Gagliardo--Nirenberg
inequality, for every $1\leq p<\infty$.

Applying it in the form \[
\displaystyle\sum_{|\beta|=j}\|D^\beta v\|_{L^p(\mathbb R^n)}
\leq
A
\|v\|_{L^p(\mathbb R^n)}^{1-\frac{j}{m}}
\left(
\displaystyle\sum_{|\gamma|=m}\|D^\gamma v\|_{L^p(\mathbb R^n)}
\right)^{\frac{j}{m}},
\qquad v\in C_c^\infty(\mathbb R^n),
\] gives, for \(v=\widetilde{F_{\alpha}^{a}\circ\phi_{\alpha}^{-1}}\), \[
\displaystyle\sum_{|\beta|=j}
\|D^\beta(F_{\alpha}^{a}\circ\phi_{\alpha}^{-1})\|_{L^p(\Omega_{\alpha})}
\leq
A
\|F_{\alpha}^{a}\circ\phi_{\alpha}^{-1}\|_{L^p(\Omega_{\alpha})}^{1-\frac{j}{m}}
\left(
\displaystyle\sum_{|\gamma|=m}
\|D^\gamma(F_{\alpha}^{a}\circ\phi_{\alpha}^{-1})\|_{L^p(\Omega_{\alpha})}
\right)^{\frac{j}{m}}.
\] Summing over \(a\) and using Hölder's inequality from Proposition~\ref{desigualdad de holder} for the finite counting measure, with exponents \(\displaystyle\frac{m}{m-j}\) and \(\displaystyle\frac{m}{j}\), gives \[
\displaystyle\sum_{a=1}^{r}\displaystyle\sum_{|\beta|=j}
\|D^\beta(F_{\alpha}^{a}\circ\phi_{\alpha}^{-1})\|_{L^p(\Omega_{\alpha})}
\leq
A
\left(
\displaystyle\sum_{a=1}^{r}\|F_{\alpha}^{a}\circ\phi_{\alpha}^{-1}\|_{L^p(\Omega_{\alpha})}
\right)^{1-\frac{j}{m}}
\left(
\displaystyle\sum_{a=1}^{r}\displaystyle\sum_{|\gamma|=m}
\|D^\gamma(F_{\alpha}^{a}\circ\phi_{\alpha}^{-1})\|_{L^p(\Omega_{\alpha})}
\right)^{\frac{j}{m}}.
\]

We now relate the Euclidean partial derivatives of the components to the covariant derivatives of \(\psi_{\alpha}\mathbf{F}\). By Lemma~\ref{lem:local-expression-higher-order}, for each \(s\in\{1,\dots,m\}\), the components of \(\nabla^s(\psi_{\alpha}\mathbf{F})\) in the chart \((U_{\alpha},\phi_{\alpha})\) and frame \((\mathbf{e}_{1,\alpha},\dots,\mathbf{e}_{r,\alpha})\) are given by \[
(\nabla^s(\psi_{\alpha}\mathbf{F}))^a_{i_1\dots i_s}
=
\partial_{i_1}\cdots\partial_{i_s}F_{\alpha}^{a}
+
\displaystyle\sum_{b=1}^{r}\displaystyle\sum_{|\beta|\leq s-1}
(B_{\alpha,\beta})^a_{b\,i_1\dots i_s}\partial^\beta F_{\alpha}^{b},
\] where the coefficients \((B_{\alpha,\beta})^a_{b\,i_1\dots i_s}\) are smooth on \(U_{\alpha}\). Since \(\supp(\psi_{\alpha}\mathbf{F})\subseteq U_{\alpha}\), all calculations take place on a compact subset of \(U_{\alpha}\), and therefore those coefficients, as well as the coefficients converting Riemannian measure to Euclidean measure in coordinates, are uniformly bounded. This triangular formula shows that there are constants \(P_{\alpha,s},Q_{\alpha,s}>0\), independent of \(\mathbf{F}\), such that \[
\|\nabla^s(\psi_{\alpha}\mathbf{F})\|_{L^p(M,T^{(0,s)}(TM)\otimes \mathbf{E})}
\leq
P_{\alpha,s}
\displaystyle\sum_{a=1}^{r}\displaystyle\sum_{|\beta|\leq s}
\|D^\beta(F_{\alpha}^{a}\circ\phi_{\alpha}^{-1})\|_{L^p(\Omega_{\alpha})},
\] and \[
\displaystyle\sum_{a=1}^{r}\displaystyle\sum_{|\beta|\leq s}
\|D^\beta(F_{\alpha}^{a}\circ\phi_{\alpha}^{-1})\|_{L^p(\Omega_{\alpha})}
\leq
Q_{\alpha,s}
\displaystyle\sum_{\ell=0}^{s}
\|\nabla^\ell(\psi_{\alpha}\mathbf{F})\|_{L^p(M,T^{(0,\ell)}(TM)\otimes \mathbf{E})}.
\] In particular, using these estimates with \(s=j\), also using the preceding Euclidean inequality for the orders \(1,\dots,j\), and noting that order \(0\) corresponds to the \(L^p(\Omega_\alpha)\) norm, we obtain a constant \(R_{\alpha}>0\), independent of \(\mathbf{F}\), such that \[
\|\nabla^j(\psi_{\alpha}\mathbf{F})\|_{L^p(M,T^{(0,j)}(TM)\otimes \mathbf{E})}
\leq
R_{\alpha}
\|\psi_{\alpha}\mathbf{F}\|_{L^p(M,\mathbf{E})}^{1-\frac{j}{m}}
\left(
\displaystyle\sum_{\ell=0}^{m}
\|\nabla^\ell(\psi_{\alpha}\mathbf{F})\|_{L^p(M,T^{(0,\ell)}(TM)\otimes \mathbf{E})}
\right)^{\frac{j}{m}}.
\] Since the cover is finite, taking \(R:=\displaystyle\max\{R_{\alpha}\mid 1\leq \alpha\leq N\}\) gives the same estimate with \(R\) in place of \(R_{\alpha}\).

We now sum over \(\alpha\). Since \(\mathbf{F}=\displaystyle\sum_{\alpha=1}^{N}\psi_{\alpha}\mathbf{F}\), the triangle inequality gives \[
\|\nabla^j\mathbf{F}\|_{L^p(M,T^{(0,j)}(TM)\otimes \mathbf{E})}
\leq
\displaystyle\sum_{\alpha=1}^{N}
\|\nabla^j(\psi_{\alpha}\mathbf{F})\|_{L^p(M,T^{(0,j)}(TM)\otimes \mathbf{E})}.
\] By the preceding estimate and Hölder's inequality from Proposition~\ref{desigualdad de holder} for the finite counting measure, \[
\|\nabla^j\mathbf{F}\|_{L^p(M,T^{(0,j)}(TM)\otimes \mathbf{E})}
\leq
R
\left(
\displaystyle\sum_{\alpha=1}^{N}\|\psi_{\alpha}\mathbf{F}\|_{L^p(M,\mathbf{E})}
\right)^{1-\frac{j}{m}}
\left(
\displaystyle\sum_{\alpha=1}^{N}
\displaystyle\sum_{\ell=0}^{m}
\|\nabla^\ell(\psi_{\alpha}\mathbf{F})\|_{L^p(M,T^{(0,\ell)}(TM)\otimes \mathbf{E})}
\right)^{\frac{j}{m}}.
\]

It remains to control the terms involving \(\psi_{\alpha}\mathbf{F}\) by \(\mathbf{F}\). For the zeroth-order term, since each \(\psi_{\alpha}\) is smooth and the family is finite, there is \(D_0>0\), independent of \(\mathbf{F}\), such that \[
\displaystyle\sum_{\alpha=1}^{N}\|\psi_{\alpha}\mathbf{F}\|_{L^p(M,\mathbf{E})}
\leq
D_0\|\mathbf{F}\|_{L^p(M,\mathbf{E})}.
\] On the other hand, for each \(s\leq m\), the Leibniz rule for the induced connection on \(T^{(0,s)}(TM)\otimes \mathbf{E}\) gives \[
|\nabla^s(\psi_{\alpha}\mathbf{F})|_{\mathbf{g},\mathbf{h}_{\mathbf{E}}}
\leq
\displaystyle\sum_{t=0}^{s}\binom{s}{t}
|\nabla^{s-t}\psi_{\alpha}|_{\mathbf{g}}\,
|\nabla^t\mathbf{F}|_{\mathbf{g},\mathbf{h}_{\mathbf{E}}}.
\] Since \(M\) is compact and the family \(\{\psi_{\alpha}\}_{\alpha=1}^{N}\) is finite, the quantities \(\||\nabla^{s-t}\psi_{\alpha}|_{\mathbf{g}}\|_{L^\infty(M)}\), with \(0\leq t\leq s\leq m\) and \(1\leq\alpha\leq N\), are uniformly bounded. Consequently, there is \(D_m>0\), independent of \(\mathbf{F}\), such that \[
\displaystyle\sum_{\alpha=1}^{N}\displaystyle\sum_{\ell=0}^{m}
\|\nabla^\ell(\psi_{\alpha}\mathbf{F})\|_{L^p(M,T^{(0,\ell)}(TM)\otimes \mathbf{E})}
\leq
D_m
\displaystyle\sum_{\ell=0}^{m}
\|\nabla^\ell \mathbf{F}\|_{L^p(M,T^{(0,\ell)}(TM)\otimes \mathbf{E})}.
\] Substituting these estimates gives \[
\|\nabla^j\mathbf{F}\|_{L^p(M,T^{(0,j)}(TM)\otimes \mathbf{E})}
\leq
C_1
\|\mathbf{F}\|_{L^p(M,\mathbf{E})}^{1-\frac{j}{m}}
\left(
\displaystyle\sum_{\ell=0}^{m}
\|\nabla^\ell \mathbf{F}\|_{L^p(M,T^{(0,\ell)}(TM)\otimes \mathbf{E})}
\right)^{\frac{j}{m}}.
\]

Finally, we use equivalence of the Sobolev norm with the norm given by the zeroth-order term and the covariant derivative of highest order; that is, there is \(C_2>0\) such that \[
\displaystyle\sum_{\ell=0}^{m}
\|\nabla^\ell \mathbf{F}\|_{L^p(M,T^{(0,\ell)}(TM)\otimes \mathbf{E})}
\leq
C_2
\left(
\|\mathbf{F}\|_{L^p(M,\mathbf{E})}
+
\|\nabla^m\mathbf{F}\|_{L^p(M,T^{(0,m)}(TM)\otimes \mathbf{E})}
\right).
\] Indeed, write $S=\sum_{\ell=0}^m\|\nabla^\ell\mathbf F\|_{L^p(M,T^{(0,\ell)}(TM)\otimes\mathbf E)}$.
The inequality just proved, applied separately to each $0<\ell<m$,
and Young's inequality give
\[
\|\nabla^\ell\mathbf F\|_{L^p(M,T^{(0,\ell)}(TM)\otimes\mathbf E)}
\leq\varepsilon S+C_{\ell,\varepsilon}\|\mathbf F\|_{L^p(M,\mathbf E)}.
\]
Sum these inequalities and choose $(m-1)\varepsilon\leq1/2$.
Moving the resulting term $S/2$ to the left proves the asserted
norm equivalence. This argument uses only the estimates already
established, and does not require a compact embedding. Therefore, \[
\|\nabla^{j}\mathbf{F}\|_{L^p(M,T^{(0,j)}(TM)\otimes \mathbf{E})}
\leq
C
\|\mathbf{F}\|_{L^p(M,\mathbf{E})}^{1-\frac{j}{m}}
\left(
\|\mathbf{F}\|_{L^p(M,\mathbf{E})}
+
\|\nabla^{m}\mathbf{F}\|_{L^p(M,T^{(0,m)}(TM)\otimes \mathbf{E})}
\right)^{\frac{j}{m}}
\] for every \(\mathbf{F}\in\Gamma(\mathbf{E})\).

Now let \(\mathbf{F}\in W^{m,p}(M,\mathbf{E})\). By Theorem~\ref{meyers-serrin-haz}, there is a sequence \((\mathbf{F}_k)_{k\in\mathbb N}\subseteq \Gamma(\mathbf{E})\cap W^{m,p}(M,\mathbf{E})\) such that \(\mathbf{F}_k\to \mathbf{F}\) in \(W^{m,p}(M,\mathbf{E})\). In particular, \(\mathbf{F}_k\to \mathbf{F}\) in \(L^p(M,\mathbf{E})\), \(\nabla^j\mathbf{F}_k\to\nabla_w^j\mathbf{F}\) in \(L^p(M,T^{(0,j)}(TM)\otimes \mathbf{E})\), and \(\nabla^m\mathbf{F}_k\to\nabla_w^m\mathbf{F}\) in \(L^p(M,T^{(0,m)}(TM)\otimes \mathbf{E})\). Applying the inequality already proved to \(\mathbf{F}_k\) gives \[
\|\nabla^{j}\mathbf{F}_k\|_{L^p(M,T^{(0,j)}(TM)\otimes \mathbf{E})}
\leq
C
\|\mathbf{F}_k\|_{L^p(M,\mathbf{E})}^{1-\frac{j}{m}}
\left(
\|\mathbf{F}_k\|_{L^p(M,\mathbf{E})}
+
\|\nabla^{m}\mathbf{F}_k\|_{L^p(M,T^{(0,m)}(TM)\otimes \mathbf{E})}
\right)^{\frac{j}{m}}.
\] Letting \(k\to\infty\), we conclude that \[
\|\nabla^{j}\mathbf{F}\|_{L^p(M,T^{(0,j)}(TM)\otimes \mathbf{E})}
\leq
C
\|\mathbf{F}\|_{L^p(M,\mathbf{E})}^{1-\frac{j}{m}}
\left(
\|\mathbf{F}\|_{L^p(M,\mathbf{E})}
+
\|\nabla_w^{m}\mathbf{F}\|_{L^p(M,T^{(0,m)}(TM)\otimes \mathbf{E})}
\right)^{\frac{j}{m}}.
\] This completes the proof. \end{proof}

\section{Sobolev spaces on manifolds with boundary} \label{sec:sobolev-haces-frontera}

In many applications, boundary conditions arise naturally when formulating a PDE. Sobolev spaces allow us to pass from a differential problem to its weak formulation; however, when the manifold has boundary, integration-by-parts formulas produce additional terms on $\partial M$. The identity with test sections supported in the interior determines the formal differential adjoint, as explained in Remark~\ref{obs:adjunto-formal-frontera-dominio}. Boundary conditions enter when choosing a domain for the realization of the operator in $L^2(M,\mathbf E)$ and determining the domain of its Hilbert space adjoint.

For this purpose, introduce the notation \[
\Gamma_{c,\operatorname{Int}(M)}(\mathbf{E})
:=
\left\{
\mathbf{u}\in\Gamma(\mathbf{E})
\middle|
\operatorname{supp}(\mathbf{u})\Subset\operatorname{Int}(M)
\right\}.
\] That is, $\Gamma_{c,\operatorname{Int}(M)}(\mathbf{E})$ denotes the space of smooth sections whose compact support lies in the interior of $M$.

Under this convention, integration-by-parts formulas produce no boundary terms when the test sections belong to $\Gamma_{c,\operatorname{Int}(M)}$. For example, in Theorem~\ref{teo:green-m-esima-covariante-E}, if we take $\mathbf{F}\in \Gamma(T^{(k,l)}(TM)\otimes \mathbf{E})$ and $\mathbf{G}\in \Gamma_{c,\operatorname{Int}(M)}(T^{(k,l+s)}(TM)\otimes \mathbf{E})$, the boundary term vanishes and the formula reduces to \[
\renewcommand{\CancelColor}{\color{red}}
\begin{aligned}
&\int_{M}\langle \nabla^{s}\mathbf{F},\mathbf{G}\rangle_{\mathbf{g},\mathbf{h}_{\mathbf{E}}}\,d\lambda_{\mathbf{g}}\\
&\quad=(-1)^{s}\int_{M}\big\langle \mathbf{F},(\operatorname{tr}_{\mathbf{g}}\circ \nabla)^{s}\mathbf{G}\big\rangle_{\mathbf{g},\mathbf{h}_{\mathbf{E}}}\,d\lambda_{\mathbf{g}}\\
&\qquad+\cancelto{\color{red}{0}}{
\displaystyle\sum_{j=0}^{s-1}(-1)^{\,s-1-j}
\int_{\partial M}
\Big\langle
\nabla^{j}\mathbf{F},\,
\iota_{\boldsymbol{\nu}}\bigl((\operatorname{tr}_{\mathbf{g}}\circ \nabla)^{\,s-1-j}\mathbf{G}\bigr)
\Big\rangle_{\mathbf{g},\mathbf{h}_{\mathbf{E}}}
\,d\lambda_{\widetilde{\mathbf{g}}}
}.
\end{aligned}
\]

Thus, when defining weak derivatives on a manifold with boundary, we work with test sections belonging to $\Gamma_{c,\operatorname{Int}(M)}$. In this way, boundary conditions are not implicitly incorporated into the weak definition; they will be introduced later through the trace operator.

We next give the definitions needed to work with these Sobolev spaces.

\begin{definition}[Interior formal adjoint]\label{def:adjunto-formal-interior}\index{interior formal adjoint} Let $(M,\mathbf{g})$ be a Riemannian manifold with nonempty boundary, and let $\mathbf{E},\mathbf{F}\to M$ be smooth vector bundles equipped with bundle metrics $\mathbf{h}_{\mathbf{E}}$ and $\mathbf{h}_{\mathbf{F}}$, respectively. Let $P\in\mathbf{PDO}(\mathbf{E},\mathbf{F})$. We say that $Q\in\mathbf{PDO}(\mathbf{F},\mathbf{E})$ is an \textit{interior formal adjoint} of $P$ if, for every $\mathbf{u}\in\Gamma_{c,\operatorname{Int}(M)}(\mathbf{E})$ and every $\mathbf{v}\in\Gamma_{c,\operatorname{Int}(M)}(\mathbf{F})$, we have \[
\int_M\langle P\mathbf{u},\mathbf{v}\rangle_{\mathbf{h}_{\mathbf{F}}}\,d\lambda_{\mathbf{g}}
=
\int_M\langle \mathbf{u},Q\mathbf{v}\rangle_{\mathbf{h}_{\mathbf{E}}}\,d\lambda_{\mathbf{g}}.
\] In this case we write $Q=P_h^{*}$; the symbol $P'$ is reserved for the complex--bilinear formal transpose. \end{definition}

The adjective ``interior'' recalls that the preceding identity is tested only against sections whose supports do not meet the boundary. Thus the operator $P_h^{*}$ describes the formal Hermitian adjoint determined by the differential part of $P$, without yet choosing boundary conditions.

\begin{theorem}[Existence and uniqueness of the interior formal adjoint]\label{teo:adjunto-formal-interior}\index{existence and uniqueness of the interior formal adjoint} Let $(M,\mathbf{g})$ be a Riemannian manifold with nonempty boundary, and let $\mathbf{E},\mathbf{F}\to M$ be smooth vector bundles with bundle metrics $\mathbf{h}_{\mathbf{E}},\mathbf{h}_{\mathbf{F}}$. If $P\in\mathbf{PDO}^{(m)}(\mathbf{E},\mathbf{F})$, then $P$ admits a unique interior formal adjoint \[
P_h^{*}\in\mathbf{PDO}^{(m)}(\mathbf{F},\mathbf{E}).
\] Moreover, let $(U,\phi)$ be a chart, and let $\mathbf{e}_1,\dots,\mathbf{e}_r$ and $\mathbf{f}_1,\dots,\mathbf{f}_s$ be local orthonormal frames of $\mathbf{E}$ and $\mathbf{F}$, respectively, defined on $U$. Set \[
D_\phi^\alpha a
:=
\bigl[D^\alpha(a\circ\phi^{-1})\bigr]\circ\phi,
\qquad a\in C^\infty(U,\mathbb K).
\] If \[
P\left(\sum_{i=1}^r u^i\mathbf e_i\right)
=
\sum_{j=1}^s
\left(
\sum_{i=1}^r\sum_{|\alpha|\leq m}
a_{\alpha i}^{\,j}D_\phi^\alpha u^i
\right)\mathbf f_j,
\] then \[
P_h^*\left(\sum_{j=1}^s v^j\mathbf f_j\right)
=
\frac{1}{\sqrt{\det(\mathbf{g})}}
\sum_{i=1}^r
\left[
\sum_{j=1}^s\sum_{|\alpha|\leq m}(-1)^{|\alpha|}
D_\phi^\alpha\!\left(
\sqrt{\det(\mathbf{g})}\,
\overline{a_{\alpha i}^{\,j}}\,v^j
\right)
\right]\mathbf e_i,
\] where $g_{ab}=\mathbf g(\boldsymbol{\partial}_a,
\boldsymbol{\partial}_b)$. In the real case the bar leaves $a_{\alpha i}^{\,j}$ unchanged. \end{theorem}

\begin{proof} Apply Theorem~\ref{teo:existencia-adjunto-hermitiano-formal} to the manifold without boundary $\operatorname{Int}(M)$ and the operator restricted to its interior. The local expression proved there depends smoothly on the coefficients of $P$, the bundle metrics, and the coordinate weight $\sqrt{\det(\mathbf{g})}$ of the Riemannian measure; since all these data are smooth up to $\partial M$, its coefficients extend smoothly to each boundary chart and define the operator displayed above throughout $M$. For test sections supported in $\operatorname{Int}(M)$, the integral identity is precisely that of the cited theorem and produces no boundary terms.

If $Q_1,Q_2$ satisfy the interior identity, the uniqueness lemma of Chapter~\ref{cap:operadores-diferenciales-parciales-haces} applied on $\operatorname{Int}(M)$ gives $Q_1=Q_2$ there. Their local coefficients are smooth up to the boundary and agree on the dense open interior; by continuity they also agree on the boundary. Therefore, $Q_1=Q_2$ on $M$. \end{proof}

\begin{definition}[Weak covariant derivative on manifolds with boundary]\label{def:derivada-covariante-debil-frontera}\index{weak covariant derivative on manifolds with boundary@weak covariant derivative on manifolds with boundary} Let $(M,\mathbf{g})$ be a Riemannian manifold with nonempty boundary and $\mathbf{E}\to M$ a smooth vector bundle equipped with a bundle metric $\mathbf{h}_{\mathbf{E}}$ and a smooth connection $\nabla^{\mathbf{E}}$. In constructing $\nabla^s$, we use the product connection induced by Levi--Civita and $\nabla^{\mathbf E}$. Let $s\geq 0$ and let $\mathbf{u}\in L^{1}_{\operatorname{loc}}(M,\mathbf{E})$. We say that a section $\mathbf{v}\in L^{1}_{\operatorname{loc}}(M,T^{(0,s)}(TM)\otimes \mathbf{E})$ is the $s$th weak covariant derivative of $\mathbf{u}$ if, for every $\boldsymbol{\sigma}\in\Gamma_{c,\operatorname{Int}(M)}(T^{(0,s)}(TM)\otimes \mathbf{E})$, we have \[
\int_M \langle \mathbf{v},\boldsymbol{\sigma}\rangle_{\mathbf{g},\mathbf{h}_{\mathbf{E}}}\,d\lambda_{\mathbf{g}}
=
\int_M \langle \mathbf{u},(\nabla^s)_h^{*}\boldsymbol{\sigma}\rangle_{\mathbf{h}_{\mathbf{E}}}\,d\lambda_{\mathbf{g}}.
\] In this case we write $\mathbf{v}=\nabla_w^s \mathbf{u}$. For $s=0$ we adopt the convention $\nabla_w^0\mathbf{u}=\mathbf{u}$. \end{definition}

In the preceding definition, $(\nabla^s)_h^*$ denotes the interior formal Hermitian adjoint of $\nabla^s$, not its bilinear transpose. This choice allows us to define weak derivatives without implicitly introducing boundary conditions.

\begin{definition}[Global and local Sobolev spaces on manifolds with boundary]\label{def:sobolev-global-local-frontera}\index{Sobolev space@Sobolev space!global and local on manifolds with boundary} Let $(M,\mathbf{g})$ be a Riemannian manifold with nonempty boundary, and let $\pi\colon \mathbf{E}\longrightarrow M$ be a smooth vector bundle equipped with a bundle metric $\mathbf{h}_{\mathbf{E}}$ and a smooth connection $\nabla^{\mathbf{E}}$. Let $m\in\mathbb N_0$ and $1\leq p\leq\infty$.

\begin{enumerate} \item \emph{Global Sobolev spaces.} Define $W^{m,p}(M,\mathbf{E})$ to be the set of sections $\mathbf{u}\in L^{1}_{\operatorname{loc}}(M,\mathbf{E})$ such that, for every $s\in\{0,\dots,m\}$, the weak covariant derivative $\nabla_w^s \mathbf{u}$ exists and belongs to $L^p(M,T^{(0,s)}(TM)\otimes \mathbf{E})$. Equip this space with the norm \[
 \|\mathbf{u}\|_{W^{m,p}(M,\mathbf{E})}
 :=
 \begin{cases}
 \displaystyle
 \left(
 \displaystyle\sum_{s=0}^{m}
 \int_M |\nabla_w^{s}\mathbf{u}|_{\mathbf{g},\mathbf{h}_{\mathbf{E}}}^{p}\,d\lambda_{\mathbf{g}}
 \right)^{\frac{1}{p}},
 & 1\leq p<\infty, \\[1.2em]
 \displaystyle
 \max_{0\leq s\leq m}
 \operatorname*{ess\,sup}_{x\in M}
 |\nabla_w^{s}\mathbf{u}(x)|_{\mathbf{g},\mathbf{h}_{\mathbf{E}}},
 & p=\infty.
 \end{cases}
 \] With this norm, $W^{m,p}(M,\mathbf{E})$ is a Banach space.

\item \emph{Local Sobolev spaces.} Define $W^{m,p}_{\operatorname{loc}}(M,\mathbf{E})$ to be the set of sections $\mathbf{u}\in L^{1}_{\operatorname{loc}}(M,\mathbf{E})$ such that, for every $s\in\{0,\dots,m\}$, the weak covariant derivative $\nabla_w^s \mathbf{u}$ exists and belongs to $L^p_{\operatorname{loc}}(M,T^{(0,s)}(TM)\otimes \mathbf{E})$.

If $K\subseteq M$ is compact, define the local seminorm \[
 \|\mathbf{u}\|_{W^{m,p}(M,\mathbf{E};K)}
 :=
 \begin{cases}
 \displaystyle
 \left(
 \displaystyle\sum_{s=0}^{m}
 \int_K |\nabla_w^{s}\mathbf{u}|_{\mathbf{g},\mathbf{h}_{\mathbf{E}}}^{p}\,d\lambda_{\mathbf{g}}
 \right)^{\frac{1}{p}},
 & 1\leq p<\infty, \\[1.2em]
 \displaystyle
 \max_{0\leq s\leq m}
 \operatorname*{ess\,sup}_{x\in K}
 |\nabla_w^{s}\mathbf{u}(x)|_{\mathbf{g},\mathbf{h}_{\mathbf{E}}},
 & p=\infty.
 \end{cases}
 \] We consider $W^{m,p}_{\operatorname{loc}}(M,\mathbf{E})$ with the family of seminorms $\{\|\cdot\|_{W^{m,p}(M,\mathbf{E};K)}\}_{K\subseteq M\ \operatorname{compact}}$. \end{enumerate} \end{definition} Restriction to the interior identifies these spaces $W^{m,p}(M,\mathbf E)$ isometrically with the spaces defined weakly on $\operatorname{Int}(M)$ using the restricted metric and connection. Indeed, the test sections are exactly those with compact support in the interior, and $\partial M$ has zero $n$-dimensional Riemannian measure, so the global norms agree. The values of a measurable section on $\partial M$ do not affect its equivalence class. In particular, Lemma~\ref{derivada debil es operador cerrado} and completeness of Sobolev spaces also apply here for every $1\leq p\leq\infty$, without a reflexivity hypothesis.

The following local characterization is proved directly using the
interior formal adjoint and the connection coefficients.

\begin{lemma}[Local characterization and norm equivalence up to the boundary]
\label{lema:meyers-serrin-haz-frontera}
Let $(M,\mathbf g)$ be a Riemannian manifold with smooth boundary,
possibly empty. Let $m\in\mathbb N_0$, $1\leq p\leq\infty$,
$(U,\phi)$ be a regular coordinate ball or half-ball, and
$\mathbf E\to M$ a real or complex bundle of rank $h$, with smooth
fiber metric and connection. The local frame and chart are assumed
defined on a neighborhood of $\overline U$. Let
$\mathbf F\in L^1_{\operatorname{loc}}(M,\mathbf E)$. Put
$B_U^+:=\phi(U\cap\operatorname{Int}(M))$ and write
$\mathbf F=\sum_{a=1}^h F^a\mathbf e_a$. Then
\[
 \mathbf F|_U\in W^{m,p}(U,\mathbf E|_U)
 \quad\Longleftrightarrow\quad
 f^a:=F^a\circ\phi^{-1}\in W^{m,p}(B_U^+)
 \quad(1\leq a\leq h).
\]
There are constants $c,C>0$, independent of $\mathbf F$, such that
\[
 c\sum_{a=1}^h\|f^a\|_{W^{m,p}(B_U^+)}
 \leq\|\mathbf F\|_{W^{m,p}(U,\mathbf E|_U)}
 \leq C\sum_{a=1}^h\|f^a\|_{W^{m,p}(B_U^+)}.
\]
\end{lemma}

\begin{proof}
Use coordinates to write the coefficients as functions on $B_U^+$.
Let
$\nabla_i^{\mathbf E}\mathbf e_b=\sum_a A_{i,b}^a\mathbf e_a$ and
$\nabla_i^{TM}\boldsymbol\partial_k
 =\sum_\ell\Gamma_{ik}^\ell\boldsymbol\partial_\ell$.
For an ordered block $I=(i_1,\ldots,i_s)$, denote by $V_{s,I}^a$
the components of $\nabla^s\mathbf F$.
With $V_0^a=f^a$, the product connection gives the recurrence
\[
 V_{s+1,iI}^a
 =\partial_iV_{s,I}^a+\sum_b A_{i,b}^aV_{s,I}^b
   -\sum_{\nu=1}^s\sum_k
       \Gamma_{ii_\nu}^kV_{s,I(\nu:k)}^a,
\]
where $I(\nu:k)$ replaces the index $i_\nu$ by $k$.
The new covariant index is placed first, in agreement with the
convention for $\nabla^{s+1}$. By induction, this recurrence defines
a matrix differential operator $P_s$ with smooth coefficients and
\[
 (P_sf)_I^a=D_If^a+
       \sum_{b=1}^h\sum_{|\beta|<s}B_{s,I,\beta,b}^aD^\beta f^b.
 \tag{*}
\]
For $s=1$, this is the connection formula. At the next step,
differentiating the principal term produces $D_iD_If^a$;
differentiating a coefficient or applying a connection coefficient
only produces derivatives of order at most $s$. This proves the
formula and smoothness of all its coefficients.

Now identify the weak operator without assuming that $f$ is smooth.
Let $\omega=\sqrt{\det(\mathbf g)}$ be the volume density, and let
$H_0,H_s$ be the positive Hermitian matrices of the metrics on
$\mathbf E$ and $T^{(0,s)}(TM)\otimes\mathbf E$ in the chosen frames.
In the real case, they are symmetric matrices.
For a matrix differential operator $P_s=\sum_\beta B_\beta D^\beta$,
its Euclidean adjoint is
\[
 P_s^\dagger\zeta
   =\sum_\beta(-1)^{|\beta|}D^\beta(B_\beta^*\zeta),
\]
where $B_\beta^*$ denotes conjugate transpose.
Integration by parts against tests compactly supported in $B_U^+$
gives the coordinate expression of the interior formal adjoint:
\[
 (P_s)_h^*\Psi
   =(\omega H_0)^{-1}P_s^\dagger(\omega H_s\Psi).
\]
These products and inverses are smooth; multiplying a compactly
supported test by any of them preserves compact support.
For every $\zeta\in C_c^\infty(B_U^+;\mathbb K^{hn^s})$,
we may therefore use $\Psi=(\omega H_s)^{-1}\zeta$.
The identity defining $\nabla_w^s\mathbf F$ is exactly equivalent to
\[
 \int_{B_U^+}\zeta^*V_s\,dx
   =\int_{B_U^+}(P_s^\dagger\zeta)^*f\,dx.
\]
This is the distributional identity $V_s=P_sf$. Conversely,
that identity for all tests $\zeta$ implies the covariant identity
by setting $\zeta=\omega H_s\Psi$.
The correspondence between tests is bijective. Their compact
supports lie in the Euclidean interior of the chart, so no
integration produces a boundary term.

If the components of $f$ belong to $W^{m,p}(B_U^+)$, formula $(*)$
shows that every $P_sf$ is an $L^p$ function representing the weak
covariant derivative. Conversely, assume that $V_s=P_sf$
belongs to $L^p$ for $0\leq s\leq m$.
At order zero we have $f\in L^p$. Once all derivatives of order
less than $s$ are known, for each block $I$ of length $s$ solve $(*)$:
\[
 D_If^a=(V_s)_I^a-
       \sum_{b=1}^h\sum_{|\beta|<s}B_{s,I,\beta,b}^aD^\beta f^b.
\]
The right-hand side belongs to $L^p$. Every multi-index of order $s$
corresponds to a block $I$, so this induction recovers all partial
derivatives through order $m$.

It remains to control the constants. On $\overline U$, all
$B_{s,I,\beta,b}^a$ are bounded, $\omega$ has positive upper and
lower bounds, and the eigenvalues of $H_s$ have uniform positive
bounds for $0\leq s\leq m$.
Formula $(*)$ and its solved version, applied inductively, give
the two estimates
\[
 \begin{aligned}
 \sum_{s=0}^m\sum_{a,I}\|(V_s)_I^a\|_{L^p(B_U^+)}
 &\leq C_1\sum_{a=1}^h\sum_{|\beta|\leq m}
                       \|D^\beta f^a\|_{L^p(B_U^+)},\\
 \sum_{a=1}^h\sum_{|\beta|\leq m}\|D^\beta f^a\|_{L^p(B_U^+)}
 &\leq C_2\sum_{s=0}^m\sum_{a,I}\|(V_s)_I^a\|_{L^p(B_U^+)}.
 \end{aligned}
\]
Here, in every sum of order $s$, $I$ ranges over the $n^s$ blocks
of that length. Equivalence of norms on these finite families,
the bounds on $H_s$, and those on $\omega$ convert these estimates
into the two inequalities in the statement. For $p=\infty$, the
same argument uses essential norms and finite maxima. The boundary
has measure zero, so the interior integrals give the norms on $U$.
\end{proof}

\begin{lemma}[Zero extension of sections with localized support]
\label{lem:frontera-extension-cero-secciones-soporte}
Let $M$ and $\mathbf E$ be as in the preceding lemma,
$U\subseteq M$ open, $K\subseteq U$ compact,
$m\in\mathbb N_0$, and $1\leq p\leq\infty$.
If $\mathbf v\in W^{m,p}(U,\mathbf E|_U)$ vanishes almost everywhere
outside $K$, its zero extension $Z_U\mathbf v$ to $M$ satisfies
\[
 \nabla_w^s(Z_U\mathbf v)=Z_U(\nabla_w^s\mathbf v),\qquad
 \|Z_U\mathbf v\|_{W^{m,p}(M,\mathbf E)}
  =\|\mathbf v\|_{W^{m,p}(U,\mathbf E|_U)},\qquad 0\leq s\leq m.
\]
The compact set $K$ may meet $\partial M$.
\end{lemma}

\begin{proof}
The local weak identity, applied to tests supported in $U\setminus K$,
shows that $\nabla_w^s\mathbf v$ vanishes almost everywhere there.
Choose $\chi\in C_c^\infty(U)$ equal to one on a neighborhood of $K$.
For a test section $\boldsymbol\sigma$ compactly supported in
$\operatorname{Int}(M)$, the product $\chi\boldsymbol\sigma|_U$
has compact support in $U\cap\operatorname{Int}(M)$ and is an
admissible local test.
Locality of the differential operator $(\nabla^s)_h^*$ implies
that $(\nabla^s)_h^*(\chi\boldsymbol\sigma)$ and
$(\nabla^s)_h^*\boldsymbol\sigma$ agree near $K$. Thus,
\[
 \begin{aligned}
 \int_M\langle Z_U(\nabla_w^s\mathbf v),\boldsymbol\sigma\rangle
                          \,d\lambda_{\mathbf g}
 &=\int_U\langle\nabla_w^s\mathbf v,
                         \chi\boldsymbol\sigma\rangle\,d\lambda_{\mathbf g}\\
 &=\int_U\langle\mathbf v,
                    (\nabla^s)_h^*(\chi\boldsymbol\sigma)\rangle
                          \,d\lambda_{\mathbf g}\\
 &=\int_M\langle Z_U\mathbf v,
                         (\nabla^s)_h^*\boldsymbol\sigma\rangle
                          \,d\lambda_{\mathbf g}.
 \end{aligned}
\]
The inner products are those induced on the corresponding bundle.
This chain proves the global weak identity. Each extension preserves
the $L^p$ norm of the derivative and hence the full norm.
The extension crosses the artificial edge of $U$ within $M$.
\end{proof}

The following lemmas distinguish the relative edge of a chart from
the faces of its model. Write
\[
 O_{n,r}:=\mathbb R^{n-r}\times(0,\infty)^r,\qquad
 W^{m,p}(G_{n,r}):=W^{m,p}(O_{n,r}),\qquad 0\leq r\leq n.
\]
For a relatively open set $\Omega\subseteq G_{n,r}$, put
$\Omega_+:=\Omega\cap O_{n,r}$. Essential support is taken in the
closed model: containment in $K\subseteq\Omega$ means that the
function vanishes almost everywhere in $\Omega_+\setminus K$.
The compact set $K$ may meet any of the faces.

\begin{lemma}[Cutoffs and zero extension within the orthant]
\label{lem:esquinas-cortes-extension-cero}
Let $\Omega\subseteq G_{n,r}$ be relatively open,
$K\subseteq\Omega$ compact, $m\in\mathbb N_0$, and
$1\leq p\leq\infty$. There is $\eta\in C_c^\infty(\mathbb R^n)$
equal to one near $K$ such that
$\operatorname{supp}\eta\cap G_{n,r}\Subset\Omega$.
If $v\in W^{m,p}(\Omega_+)$ vanishes almost everywhere outside $K$,
its zero extension $Z_\Omega v$ to $O_{n,r}$ belongs to
$W^{m,p}(O_{n,r})$, and
\[
 D^\beta Z_\Omega v=Z_\Omega(D^\beta v),\qquad
 \|Z_\Omega v\|_{W^{m,p}(O_{n,r})}
   =\|v\|_{W^{m,p}(\Omega_+)},\qquad |\beta|\leq m.
\]
Moreover, multiplication by a smooth function with bounded derivatives
through order $m$ is continuous on each of these spaces.
\end{lemma}

\begin{proof}
If $K$ is empty, the assertions concerning $v$ are immediate.
Otherwise $G_{n,r}\setminus\Omega$ is closed in $\mathbb R^n$ and
disjoint from $K$. Choose a bounded Euclidean neighborhood of $K$
whose closure avoids this closed set, and a smooth cutoff equal to
one near $K$ and supported in that neighborhood. Its intersection
with $G_{n,r}$ is a compact subset of $\Omega$.

The weak Leibniz rule follows first at order one. For a test function
$\varphi$ compactly supported in the open set under consideration,
\[
 \begin{aligned}
 \int_{\Omega_+} \eta v\,\partial_i\varphi\,dx
 &=\int_{\Omega_+}v\,\partial_i(\eta\varphi)\,dx
       -\int_{\Omega_+}v\,(\partial_i\eta)\varphi\,dx\\
 &=-\int_{\Omega_+}
       \bigl(\eta\,\partial_i v+(\partial_i\eta)v\bigr)\varphi\,dx.
 \end{aligned}
\]
All integrals are over the same Euclidean open set.
Iterating this identity and applying Pascal's rule to the binomial
coefficients gives
\[
 D^\beta(\eta v)=
 \sum_{\gamma\leq\beta}\binom{\beta}{\gamma}
       D^\gamma\eta\,D^{\beta-\gamma}v.
\]
If $A:=\max_{|\gamma|\leq m}
\|D^\gamma\eta\|_{L^\infty(\Omega_+)}$, the triangle inequality
for each derivative and equivalence of norms on finite families give
\[
 \|\eta v\|_{W^{m,p}(\Omega_+)}
 \leq C_{n,m,p}A\|v\|_{W^{m,p}(\Omega_+)}.
\]
For $p=\infty$, use maxima instead of sums of powers.

Locality of weak derivatives implies that $D^\beta v=0$ almost
everywhere in $\Omega_+\setminus K$. Indeed, for
$\varphi\in C_c^\infty(\Omega_+\setminus K)$,
\[
 \int_{\Omega_+\setminus K}D^\beta v\,\varphi\,dx
 =(-1)^{|\beta|}
   \int_{\Omega_+\setminus K}v\,D^\beta\varphi\,dx=0;
\]
uniqueness of the distribution represented by a locally integrable
function proves the assertion.

Now let $\varphi\in C_c^\infty(O_{n,r})$.
The function $(\eta\varphi)|_{\Omega_+}$ belongs to
$C_c^\infty(\Omega_+)$: the cutoff stays away from the relative edge
of $\Omega$, and the test function stays away from every face.
Since $\eta=1$ near $K$,
\[
 \begin{aligned}
 \int_{O_{n,r}}Z_\Omega v\,D^\beta\varphi\,dx
 &=\int_{\Omega_+}v\,D^\beta(\eta\varphi)\,dx\\
 &=(-1)^{|\beta|}
      \int_{\Omega_+}D^\beta v\,\eta\varphi\,dx\\
 &=(-1)^{|\beta|}
      \int_{O_{n,r}}Z_\Omega(D^\beta v)\,\varphi\,dx.
 \end{aligned}
\]
This identifies each weak derivative. Zero extension preserves its
$L^p$ norm, including $p=\infty$; combining the derivative orders
gives the stated norm equality.
\end{proof}

\begin{lemma}[Smooth approximation up to the faces with controlled support]
\label{lem:esquinas-aproximacion-ortante-soporte}
Let $m\in\mathbb N_0$ and $1\leq p<\infty$. The restrictions of
$C_c^\infty(\mathbb R^n)$ to $O_{n,r}$ are dense in
$W^{m,p}(O_{n,r})$. If $\Omega$ and $K$ are as in the preceding
lemma and $v\in W^{m,p}(\Omega_+)$ vanishes outside $K$, there are
$f_k\in C_c^\infty(\mathbb R^n)$ and a fixed compact set
$L\subseteq\Omega$ such that
\[
 \operatorname{supp}f_k\cap G_{n,r}\subseteq L,\qquad
 \|f_k|_{\Omega_+}-v\|_{W^{m,p}(\Omega_+)}\longrightarrow0.
\]
Their restrictions are smooth on all of $\Omega$, including its faces.
\end{lemma}

\begin{proof}
Let $u\in W^{m,p}(O_{n,r})$. Denote by $u^0$ and
$(D^\beta u)^0$ the zero extensions as functions in
$L^p(\mathbb R^n)$. Let
$e=(0,\ldots,0,1,\ldots,1)$, with $r$ final ones, and let
$\rho\in C_c^\infty(B_{\mathrm{euc}}(0,1))$ be a nonnegative
mollifier of integral one. Define on all of $\mathbb R^n$
\[
 u_\delta(x):=(\rho_{\delta/4}*u^0)(x+\delta e),\qquad \delta>0.
\]
This function is smooth. If $x\in G_{n,r}$, all points $y$ in the
support of $y\mapsto\rho_{\delta/4}(x+\delta e-y)$ lie in
$O_{n,r}$, because their normal coordinates exceed $3\delta/4$.
That function is an admissible test for the weak derivatives of $u$.
Differentiating the convolution and applying the weak identity to
this test gives, for $|\beta|\leq m$,
\[
 D^\beta u_\delta(x)
   =(\rho_{\delta/4}*(D^\beta u)^0)(x+\delta e),\qquad x\in G_{n,r}.
\]
The zero extensions have been used here as $L^p$ functions;
the derivative identity was computed inside the open orthant.
Translation invariance of the norm, mollifier approximation, and
continuity of translations in $L^p$ give
\[
 \begin{aligned}
 \|D^\beta u_\delta-D^\beta u\|_{L^p(O_{n,r})}
 &\leq
 \|\rho_{\delta/4}*(D^\beta u)^0-(D^\beta u)^0\|_{L^p(\mathbb R^n)}\\
 &\quad+
 \|(D^\beta u)^0(\,\cdot+\delta e)-(D^\beta u)^0\|_{L^p(\mathbb R^n)}
 \longrightarrow0.
 \end{aligned}
\]
There are finitely many multi-indices of order at most $m$, so
$u_\delta|_{O_{n,r}}\to u$ in $W^{m,p}(O_{n,r})$.

To obtain compact support, choose a cutoff $\chi$ equal to one on
the unit ball and supported in the ball of radius two, and put
$\chi_R(x)=\chi(x/R)$. In the Leibniz formula for
$D^\beta(\chi_Ru-u)$, the term $(\chi_R-1)D^\beta u$ tends to zero
by dominated convergence. Every term with $|\gamma|\geq1$
derivatives on the cutoff satisfies
\[
 \|D^\gamma\chi_R\,D^{\beta-\gamma}u\|_{L^p(O_{n,r})}
 \leq R^{-|\gamma|}\|D^\gamma\chi\|_{L^\infty(\mathbb R^n)}
                    \|D^{\beta-\gamma}u\|_{L^p(O_{n,r})}\longrightarrow0.
\]
Thus $\chi_Ru\to u$ in the Sobolev norm. The measurable zero extension
of $\chi_Ru$ has compact support, so its translated convolution
belongs to $C_c^\infty(\mathbb R^n)$.
Choose $R_k$ so that the first error is smaller than $1/(2k)$,
and then $\delta_k$ so that the second is smaller than $1/(2k)$.
This proves density.

For the localized assertion, put $w:=Z_\Omega v$ and choose the
cutoff $\eta$ from Lemma~\ref{lem:esquinas-cortes-extension-cero}.
Since $w$ is supported in $K$, its translated convolutions
$w_{\delta_k}$ have compact support. Define
$f_k:=\eta w_{\delta_k}$ and
$L:=\operatorname{supp}\eta\cap G_{n,r}$. The identity $\eta w=w$
and the multiplier estimate of the same lemma imply
\[
 \|f_k|_{\Omega_+}-v\|_{W^{m,p}(\Omega_+)}
 \leq C_\eta\|w_{\delta_k}|_{O_{n,r}}-w\|_{W^{m,p}(O_{n,r})}
 \longrightarrow0.
\]
The compact set $L$ is independent of $k$ and lies in $\Omega$.
Each $f_k$ is smooth on a Euclidean open set containing all the
faces, which proves the final assertion.
\end{proof}

\begin{lemma}[Sobolev extension from an orthant]
\label{lem:extension-sobolev-ortante-esquinas}
Let $m\in\mathbb N_0$, $0\leq r\leq n$, and
$1\leq p\leq\infty$. There is a continuous linear operator
\[
 \mathcal E_{m,r}\colon W^{m,p}(G_{n,r})\longrightarrow
 W^{m,p}(\mathbb R^n),\qquad
 (\mathcal E_{m,r}u)|_{O_{n,r}}=u.
\]
For $m\geq1$, it can be chosen using finitely many reflections,
with coefficients depending only on $m$, and
\[
 \|\mathcal E_{m,r}u\|_{W^{m,p}(\mathbb R^n)}
 \leq C_{m,p}^{\,r}\|u\|_{W^{m,p}(O_{n,r})}.
\]
If $u$ is supported in a fixed compact set $K\subseteq G_{n,r}$,
its extension is supported in a fixed Euclidean compact set
$\widehat K$, independent of $u$.
\end{lemma}

\begin{proof}
For $r=0$, take the identity. For $m=0$, take zero extension, of
norm one. Assume $r\geq1$ and $m\geq1$.
The system
\[
 \sum_{\ell=1}^m a_\ell(-\ell)^b=1,\qquad 0\leq b\leq m-1,
\]
has a unique solution, because its Vandermonde determinant is a
product of differences of the distinct nodes $-1,\ldots,-m$.
Write $x=(y,t)$ and
$D:=\mathbb R^{n-r}\times(0,\infty)^{r-1}$. For a restriction $u$
of a function in $C_c^\infty(\mathbb R^n)$, define
\[
 (Eu)(y,t)=
 \begin{cases}
 u(y,t),&t>0,\\
 \displaystyle\sum_{\ell=1}^m a_\ell u(y,-\ell t),&t<0.
 \end{cases}
\]
For $\alpha=(\alpha',b)$ of order at most $m-1$, the limits of
$D^\alpha Eu$ at $t=0$ agree, because
\[
 \sum_{\ell=1}^m a_\ell(-\ell)^b
             D^{\alpha'}\partial_t^b u(y,0)
       =D^{\alpha'}\partial_t^b u(y,0).
\]
If $h$ is a derivative of order at most $m-1$ and
$\varphi\in C_c^\infty(D\times\mathbb R)$, Fubini and
one-dimensional integration by parts on the two sides give
\[
 \begin{aligned}
 \int_{D\times(\mathbb R\setminus\{0\})}h\,\partial_t\varphi\,dy\,dt
 &=-\int_{D\times(\mathbb R\setminus\{0\})}
                      (\partial_t h)\varphi\,dy\,dt\\
 &\quad+\int_D
     \bigl(h(y,0^-)-h(y,0^+)\bigr)\varphi(y,0)\,dy.
 \end{aligned}
\]
The face integral vanishes because the limits agree.
For a derivative in $y$, the test is supported inside $D$, and
integration by parts there produces no boundary integral.
Starting with $h=Eu$ and applying these identities successively
gives all weak derivatives of order $|\alpha|\leq m$:
\[
 D^\alpha Eu(y,t)=
 \begin{cases}
 D^\alpha u(y,t),&t>0,\\
 \displaystyle\sum_{\ell=1}^m a_\ell(-\ell)^b
                           D^\alpha u(y,-\ell t),&t<0.
 \end{cases}
\]
At every step, the derivative being integrated has order at most
$m-1$, exactly the range in which matching was verified.

For $p<\infty$, the substitution $s=-\ell t$ and Minkowski give
\[
 \|D^\alpha Eu\|_{L^p(D\times(-\infty,0))}
 \leq\sum_{\ell=1}^m |a_\ell|\ell^{b-\frac{1}{p}}
                     \|D^\alpha u\|_{L^p(D\times(0,\infty))}.
\]
The norm over the union of both sides is bounded by the sum of their
norms. Combining the multi-indices, we may take
$C_{m,p}=1+\sum_{\ell=1}^m|a_\ell|\ell^{m-\displaystyle\frac{1}{p}}$.
Lemma~\ref{lem:esquinas-aproximacion-ortante-soporte} gives a smooth
sequence converging to any $u\in W^{m,p}(O_{n,r})$.
The estimate for differences makes its extensions Cauchy in
$W^{m,p}(D\times\mathbb R)$; completeness defines $Eu$.
The reflection formula is also continuous in $L^p$, by changes
of variables, and hence agrees almost everywhere with this limit.
Restriction to the positive side still gives $u$.
Repeat the construction in each remaining normal coordinate:
every step removes one positivity restriction and has the same
bound. The composition is $\mathcal E_{m,r}$, with norm at most
$C_{m,p}^{\,r}$.

For $p=\infty$, use the finite formula directly. On any compact
subset of $D\times\mathbb R$, choose an orthant cutoff equal to
one near its preimages under all the reflections. Its product
with $u$ belongs to $W^{m,1}(O_{n,r})$: the derivatives are finite
sums of bounded products supported in a set of finite measure.
The construction for $p=1$ identifies the derivatives of its
extension on the compact set. There the cutoff and its positive
order derivatives are respectively one and zero on the preimages,
so the resulting formulas are those for $Eu$.
Since every test has compact support, these identities give the
weak derivatives of $Eu$ on all of $D\times\mathbb R$.
The essential norm estimate uses $\ell^b$ instead of
$\ell^{b-1/p}$; take
$C_{m,\infty}=1+\sum_{\ell=1}^m|a_\ell|\ell^m$.
Iteration in the remaining coordinates proves the infinite case.

Finally, if
$K\subseteq[-R,R]^{n-r}\times[0,R]^r$, every nonzero value of
$u(y,-\ell t)$ with $t<0$ satisfies $|t|\leq R/\ell$.
All reflections and their compositions are supported in
$\widehat K=[-R,R]^n$. The almost everywhere formula preserves
this property for all exponents.
\end{proof}

\begin{lemma}[Localized embeddings on relative orthant domains]
\label{lem:esquinas-encajes-localizados-ortante}
Let $\Omega\subseteq G_{n,r}$ be relatively open,
$K\subseteq\Omega$ compact, $\ell>j\geq0$ integers, and
$1\leq p,q<\infty$. If
\[
 \ell-\frac{n}{p}\geq j-\frac{n}{q},
\]
there is $C>0$ such that every $v\in W^{\ell,p}(\Omega_+)$
supported in $K$ belongs to $W^{j,q}(\Omega_+)$ and satisfies
\[
 \|v\|_{W^{j,q}(\Omega_+)}
 \leq C\|v\|_{W^{\ell,p}(\Omega_+)}.
\]
Under the strict inequality, every bounded sequence of such
functions, with the same $K$, has a subsequence converging in
$W^{j,q}(\Omega_+)$.
\end{lemma}

\begin{proof}
The two extension lemmas give
$w:=Z_\Omega v\in W^{\ell,p}(O_{n,r})$ and
$z:=\mathcal E_{\ell,r}w\in W^{\ell,p}(\mathbb R^n)$, with
\[
 \|z\|_{W^{\ell,p}(\mathbb R^n)}
 \leq C_{\ell,p}^{\,r}\|v\|_{W^{\ell,p}(\Omega_+)}.
\]
All functions $z$ are supported in the same compact set $\widehat K$.
Choose a bounded open ball $B$ with $\widehat K\Subset B$.
If $\varepsilon<\operatorname{dist}(\widehat K,\partial B)$,
the convolution $\rho_\varepsilon*z$ has compact support in $B$.
Its derivatives through order $\ell$ converge to those of $z$
in $L^p(\mathbb R^n)$. Thus $z|_B\in W_0^{\ell,p}(B)$.
Theorem~\ref{prop:encaje-continuo-indices-finitos} implies
\[
 \|z|_B\|_{W^{j,q}(B)}
 \leq C_B\|z|_B\|_{W^{\ell,p}(B)}.
\]
Zero extension from $W_0^{j,q}(B)$ preserves its derivatives
and norm, by Proposition~\ref{prop:extension-cero-W0-orden-entero}.
Its restriction to $\Omega_+$ is $v$, because $z=w$ on $O_{n,r}$.
Restriction does not increase derivative norms. The preceding
two bounds prove the stated estimate.

Under the strict inequality, for a bounded sequence $(v_k)$ use
the same operators and the same ball. The sequence $(z_k|_B)$
is bounded in $W_0^{\ell,p}(B)$. Theorem~\ref{rellich kondrashov generalizado}
gives a subsequence converging to $z_\infty$ in $W_0^{j,q}(B)$.
Zero extension and restriction to $\Omega_+$ are continuous, so
\[
 \|v_k-(Z_Bz_\infty)|_{\Omega_+}\|_{W^{j,q}(\Omega_+)}
 \leq\|z_k|_B-z_\infty\|_{W^{j,q}(B)}\longrightarrow0.
\]
This proves compactness with fixed support, including all the faces.
\end{proof}

\begin{theorem}[Meyers--Serrin on open sets with uniformly Lipschitz boundary]\label{teo:meyers-serrin-euclidiano-frontera-continua}\index{Meyers Serrin@Meyers--Serrin!on Euclidean open sets with uniformly Lipschitz boundary} Let $\Omega\subseteq\mathbb R^n$ be an open set with uniformly Lipschitz continuous boundary in the sense of Definition~\ref{def:frontera-uniformemente-lipschitz}. Let $m\in\mathbb N$ and $1\leq p<\infty$. Then the restrictions to $\Omega$ of functions in $C_c^\infty(\mathbb R^n)$ form a dense subspace of $W^{m,p}(\Omega)$. \end{theorem}

\begin{proof} Let $u\in W^{m,p}(\Omega)$. By Theorem~\ref{teo:extension-sobolev-uniforme-lipschitz}, there is $v:=\mathcal Eu\in W^{m,p}(\mathbb R^n)$ such that $v\restriction_\Omega=u$ almost everywhere. We prove explicitly that $C_c^\infty(\mathbb R^n)$ is dense in $W^{m,p}(\mathbb R^n)$.

Choose $\chi\in C_c^\infty(\mathbb R^n)$ with $0\leq\chi\leq1$, $\chi=1$ on $B_{\mathrm{euc}}(0,1)$, and $\operatorname{supp}\chi\subset B_{\mathrm{euc}}(0,2)$, and set $\chi_R(x):=\chi(x/R)$ and $v_R:=\chi_Rv$. For $|\alpha|\leq m$, the weak Leibniz rule gives \[
D^\alpha v_R
=
\sum_{\beta\leq\alpha}
{\alpha\choose\beta}
D^\beta\chi_R\,D^{\alpha-\beta}v.
\] The term with $\beta=0$ converges to $D^\alpha v$ in $L^p(\mathbb R^n)$ by dominated convergence. If $\beta\neq0$, then \[
\|D^\beta\chi_R\,D^{\alpha-\beta}v\|_{L^p(\mathbb R^n)}
\leq
R^{-|\beta|}\|D^\beta\chi\|_{L^\infty(\mathbb R^n)}
\|D^{\alpha-\beta}v\|_{L^p(\mathbb R^n)}
\longrightarrow0.
\] Since the sum has finitely many terms, $v_R\to v$ in $W^{m,p}(\mathbb R^n)$.

Now fix $R$. If $\rho_\varepsilon$ is the mollifier from Theorem~\ref{teo:molificadores-euclidianos}, then $v_{R,\varepsilon}:=\rho_\varepsilon*v_R$ belongs to $C_c^\infty(\mathbb R^n)$ and \[
D^\alpha v_{R,\varepsilon}
=
\rho_\varepsilon*D^\alpha v_R,
\qquad |\alpha|\leq m.
\] Continuity of translations in $L^p(\mathbb R^n)$ and Minkowski's integral inequality imply $\rho_\varepsilon*D^\alpha v_R\to D^\alpha v_R$ in $L^p(\mathbb R^n)$ for each $|\alpha|\leq m$. Therefore, $v_{R,\varepsilon}\to v_R$ in $W^{m,p}(\mathbb R^n)$. A diagonal choice gives $(\varphi_j)\subset C_c^\infty(\mathbb R^n)$ with $\varphi_j\to v$ in $W^{m,p}(\mathbb R^n)$. Since restriction to $\Omega$ does not increase any derivative norm, \[
\|\varphi_j\restriction_\Omega-u\|_{W^{m,p}(\Omega)}
\leq
\|\varphi_j-v\|_{W^{m,p}(\mathbb R^n)}
\longrightarrow0.
\] \end{proof}

To discuss boundary conditions, we need to define an operator known as the \textit{trace operator}. This operator is first constructed on smooth sections and then extended by density. For this reason, before introducing it, we must verify that the Meyers--Serrin theorem remains valid on manifolds with boundary.

The Euclidean Meyers--Serrin theorem for open sets with uniformly Lipschitz boundary, Theorem~\ref{teo:meyers-serrin-euclidiano-frontera-continua}, allows approximation in boundary charts by restrictions of smooth functions defined on $\mathbb R^n$.

\begin{definition}\label{def:Hmp-frontera}\index{Sobolev space@Sobolev space!definition using \(C^{m,p}(M,\mathbf{E})\) on manifolds with boundary} Let $(M,\mathbf{g})$ be a Riemannian manifold with nonempty boundary and $\mathbf{E}\longrightarrow M$ a smooth vector bundle equipped with a bundle metric $\mathbf{h}_{\mathbf{E}}$ and a smooth connection $\nabla^{\mathbf{E}}$. Let $m\in\mathbb N_0$ and let $1\leq p<\infty$. Define \[
C^{m,p}(M,\mathbf{E})
:=
\left\{
\mathbf{u}\in\Gamma(\mathbf{E})
\middle|
\|\mathbf{u}\|_{W^{m,p}(M,\mathbf{E})}<\infty
\right\},
\] where, for a smooth section $\mathbf{u}$, the Sobolev norm is given by \[
\|\mathbf{u}\|_{W^{m,p}(M,\mathbf{E})}
:=
\left(
\displaystyle\sum_{s=0}^{m}
\int_M |\nabla^s \mathbf{u}|_{\mathbf{g},\mathbf{h}_{\mathbf{E}}}^{p}\,d\lambda_{\mathbf{g}}
\right)^{\frac{1}{p}}.
\] Define \[
H^{m,p}(M,\mathbf{E})
:=
\overline{C^{m,p}(M,\mathbf{E})}^{\,\|\cdot\|_{W^{m,p}(M,\mathbf{E})}}.
\] That is, $H^{m,p}(M,\mathbf{E})$ is the completion of the smooth sections with finite Sobolev norm. \end{definition}

As in the case without boundary, this definition agrees with the weak definition of Sobolev spaces. The technical difference is that, when working locally near $\partial M$, we use regular coordinate half-balls and pass to the Euclidean open set $\phi(U\cap\operatorname{Int}(M))$. \begin{theorem}[Meyers--Serrin for vector bundles on manifolds with boundary]
\label{teo:meyers-serrin-haz-frontera}
\index{meyers serrin para haces vectoriales en variedades con frontera@Meyers--Serrin for vector bundles on manifolds with boundary}
Let $(M,\mathbf g)$ be a Riemannian manifold with nonempty smooth
boundary, and let $\mathbf E\to M$ be a finite rank real or complex
vector bundle with smooth fiber metric and compatible connection.
For $m\in\mathbb N_0$ and $1\leq p<\infty$,
\[
 H^{m,p}(M,\mathbf E)=W^{m,p}(M,\mathbf E).
\]
The identification is isometric. Neither compactness nor
completeness of $M$ is required.
\end{theorem}

\begin{proof}
Smooth sections of finite norm belong to the weak space because
the tests are supported in $\operatorname{Int}(M)$.
Completeness of the weak space canonically identifies their
completion with a closed subspace of $W^{m,p}(M,\mathbf E)$.
This gives the inclusion of $H^{m,p}$ into $W^{m,p}$ and its
isometric character. We prove density for the reverse inclusion.

Put $n:=\dim M$ and $h:=\operatorname{rank}\mathbf E$.
Proposition~\ref{prop:esquinas-cartas-regulares-cubierta}, applied
to depth zero and one models, supplies regular charts
$(U_\alpha,\phi_\alpha)$, frames on neighborhoods of their closures,
and a partition $(\psi_\alpha)_{\alpha\in A}$ with compact support
in each $U_\alpha$. The set $A$ is a subset of $\mathbb N$, and
the family $(U_\alpha)$ is locally finite.
Put
\[
 \Omega_\alpha:=\phi_\alpha(U_\alpha),\qquad
 \Omega_{\alpha,+}:=\phi_\alpha(U_\alpha\cap\operatorname{Int}(M)),
 \qquad K_\alpha:=\phi_\alpha(\operatorname{supp}\psi_\alpha).
\]
The model is $G_{n,0}$ for an interior chart and $G_{n,1}$ for a
boundary chart. In both cases, $K_\alpha$ is a compact subset
of the relatively open set $\Omega_\alpha$.

Fix $\mathbf F\in W^{m,p}(M,\mathbf E)$ and write
$\mathbf F=\sum_a F_\alpha^a\mathbf e_{\alpha,a}$ on $U_\alpha$.
Lemma~\ref{lema:meyers-serrin-haz-frontera} shows that
$F_\alpha^a\circ\phi_\alpha^{-1}\in W^{m,p}(\Omega_{\alpha,+})$.
The derivatives of the cutoff $\psi_\alpha\circ\phi_\alpha^{-1}$
through order $m$ are bounded. The weak Leibniz rule in
Lemma~\ref{lem:esquinas-cortes-extension-cero} implies that
\[
 f_\alpha^a:=(\psi_\alpha F_\alpha^a)\circ\phi_\alpha^{-1}
       \in W^{m,p}(\Omega_{\alpha,+}).
\]
These functions vanish almost everywhere outside $K_\alpha$.
The reverse local characterization also shows that
$(\psi_\alpha\mathbf F)|_{U_\alpha}$ belongs to the corresponding
covariant space.

Lemma~\ref{lem:esquinas-aproximacion-ortante-soporte}, with
$r=0$ or $r=1$, supplies functions
$f_{\alpha,k}^a\in C_c^\infty(\mathbb R^n)$ such that
\[
 \sum_{a=1}^h
 \|f_{\alpha,k}^a|_{\Omega_{\alpha,+}}-f_\alpha^a
                   \|_{W^{m,p}(\Omega_{\alpha,+})}\longrightarrow0.
\]
The supports of all their restrictions to the model lie in the
same compact set $L_\alpha\subseteq\Omega_\alpha$.
Define on $U_\alpha$
\[
 \mathbf G_{\alpha,k}
    :=\sum_{a=1}^h(f_{\alpha,k}^a\circ\phi_\alpha)\mathbf e_{\alpha,a}.
\]
The Euclidean expressions are smooth up to the face, and the support
lies in $\phi_\alpha^{-1}(L_\alpha)\Subset U_\alpha$.
Consequently, extending $\mathbf G_{\alpha,k}$ by zero outside
$U_\alpha$ gives a globally smooth section of compact support.
Local norm equivalence and
Lemma~\ref{lem:frontera-extension-cero-secciones-soporte} give
\[
 \begin{aligned}
 \|\mathbf G_{\alpha,k}-\psi_\alpha\mathbf F\|_{W^{m,p}(M,\mathbf E)}
 &=\|\mathbf G_{\alpha,k}-(\psi_\alpha\mathbf F)|_{U_\alpha}
                       \|_{W^{m,p}(U_\alpha,\mathbf E|_{U_\alpha})}\\
 &\leq C_\alpha\sum_{a=1}^h
       \|f_{\alpha,k}^a|_{\Omega_{\alpha,+}}-f_\alpha^a
                       \|_{W^{m,p}(\Omega_{\alpha,+})}
 \longrightarrow0.
 \end{aligned}
\]
The local error is supported in
$\phi_\alpha^{-1}(L_\alpha)\cup\operatorname{supp}\psi_\alpha$,
a compact subset of the domain, as required by the extension lemma.

For $q\in\mathbb N$, choose an index $k_{\alpha,q}$ for each
$\alpha$ such that
\[
 \mathbf e_{\alpha,q}:=\psi_\alpha\mathbf F-\mathbf G_{\alpha,k_{\alpha,q}},
 \qquad
 \|\mathbf e_{\alpha,q}\|_{W^{m,p}(M,\mathbf E)}<2^{-(\alpha+q)}.
\]
The local constants may depend on $\alpha$; each index is chosen
after the chart has been fixed.
The sum
$\mathbf G_q:=\sum_{\alpha\in A}\mathbf G_{\alpha,k_{\alpha,q}}$
is locally finite, so it is a globally smooth section.

To justify the estimate for the sum of errors, consider
$S_N:=\sum_{\alpha\in A,\ \alpha\leq N}\mathbf e_{\alpha,q}$.
For $N>L$,
\[
 \|S_N-S_L\|_{W^{m,p}(M,\mathbf E)}
 \leq\sum_{\substack{\alpha\in A\\L<\alpha\leq N}}
                          2^{-(\alpha+q)}.
\]
The tails tend to zero, and completeness gives a limit
$S\in W^{m,p}(M,\mathbf E)$. On a neighborhood meeting only finitely
many $U_\alpha$, the restrictions of $S_N$, for sufficiently large
$N$, equal the finite sum $\mathbf F-\mathbf G_q$.
Convergence in $L^p$ identifies the limit $S$ there with that
difference. A countable cover by such neighborhoods identifies the
two on all of $M$, up to a null set.
Passing to the limit in the triangle inequality gives
\[
 \|\mathbf F-\mathbf G_q\|_{W^{m,p}(M,\mathbf E)}
 \leq\sum_{\alpha\in A}2^{-(\alpha+q)}
 \leq2^{-q}.
\]
In particular, $\mathbf G_q=\mathbf F-S$ has finite norm, with
\[
 \|\mathbf G_q\|_{W^{m,p}(M,\mathbf E)}
 \leq\|\mathbf F\|_{W^{m,p}(M,\mathbf E)}+2^{-q}.
\]
Thus $\mathbf G_q\in\Gamma(M,\mathbf E)\cap W^{m,p}(M,\mathbf E)$
and converges to $\mathbf F$ in the Sobolev norm. This proves the
reverse inclusion and the isometric equality.
\end{proof} \begin{lemma}[Continuity of multiplication by smooth functions]\label{lem:multiplicacion-funcion-suave-sobolev}\index{continuity of multiplication by smooth functions@continuity of multiplication by smooth functions} Let $(M,\mathbf{g})$ be a Riemannian manifold with or without boundary, and let $\mathbf{E}\to M$ be a smooth vector bundle equipped with a bundle metric $\mathbf{h}_{\mathbf{E}}$ and a compatible connection $\nabla^{\mathbf{E}}$. Let $m\in\mathbb N$ and $1\leq p<\infty$. Suppose that $\psi\in C^\infty(M)$ satisfies \[
\|\nabla^q\psi\|_{L^\infty(M,T^{(0,q)}(TM))}<\infty
\qquad\text{for every }0\leq q\leq m.
\] Then the operator \[
\mathbf{u}\longmapsto \psi \mathbf{u}
\] is linear and continuous from $W^{m,p}(M,\mathbf{E})$ to itself. In particular, there is a constant $C_{\psi,m}>0$ such that \[
\|\psi \mathbf{u}\|_{W^{m,p}(M,\mathbf{E})}
\leq
C_{\psi,m}\|\mathbf{u}\|_{W^{m,p}(M,\mathbf{E})}
\] for every $\mathbf{u}\in W^{m,p}(M,\mathbf{E})$.

In particular, if $M$ is compact, every function $\psi\in C^\infty(M)$ defines a continuous multiplication operator on $W^{m,p}(M,\mathbf{E})$. More generally, on a noncompact manifold it suffices to take $\psi\in C_c^\infty(M)$. \end{lemma}

\begin{proof} We first prove the estimate for smooth sections. Let $\mathbf{u}\in \Gamma(\mathbf{E})\cap W^{m,p}(M,\mathbf{E})$. For each $0\leq s\leq m$, the iterated Leibniz rule for the induced connection shows that $\nabla^s(\psi \mathbf{u})$ is a finite sum of terms constructed from the covariant derivatives of $\psi$ and $\mathbf{u}$, whose orders sum to $s$. Consequently, there is a constant $A_s>0$, depending only on $s$ and the tensor conventions, such that \[
|\nabla^s(\psi \mathbf{u})|_{\mathbf{g},\mathbf{h}_{\mathbf{E}}}
\leq
A_s\displaystyle\sum_{t=0}^{s}
|\nabla^{s-t}\psi|_{\mathbf{g}}\,|\nabla^t \mathbf{u}|_{\mathbf{g},\mathbf{h}_{\mathbf{E}}}.
\] Using the boundedness hypothesis on the derivatives of $\psi$ gives \[
|\nabla^s(\psi \mathbf{u})|_{\mathbf{g},\mathbf{h}_{\mathbf{E}}}
\leq
A_s\displaystyle\sum_{t=0}^{s}
\|\nabla^{s-t}\psi\|_{L^\infty(M,T^{(0,s-t)}(TM))}
|\nabla^t \mathbf{u}|_{\mathbf{g},\mathbf{h}_{\mathbf{E}}}.
\] Taking the $L^p\bigl(M,T^{(0,s)}(TM)\otimes \mathbf{E}\bigr)$ norm and using the triangle inequality, \[
\|\nabla^s(\psi \mathbf{u})\|_{L^p(M,T^{(0,s)}(TM)\otimes \mathbf{E})}
\leq
A_s\displaystyle\sum_{t=0}^{s}
\|\nabla^{s-t}\psi\|_{L^\infty(M,T^{(0,s-t)}(TM))}
\|\nabla^t \mathbf{u}\|_{L^p(M,T^{(0,t)}(TM)\otimes \mathbf{E})}.
\] Summing over $s\in\{0,\dots,m\}$ gives a constant $C_{\psi,m}>0$ such that \[
\|\psi \mathbf{u}\|_{W^{m,p}(M,\mathbf{E})}
\leq
C_{\psi,m}\|\mathbf{u}\|_{W^{m,p}(M,\mathbf{E})}
\] for every $\mathbf{u}\in \Gamma(\mathbf{E})\cap W^{m,p}(M,\mathbf{E})$.

Now take $\mathbf u\in W^{m,p}(M,\mathbf E)$. If
$\partial M=\varnothing$, apply Theorem~\ref{meyers-serrin-haz};
if $\partial M\neq\varnothing$, apply
Theorem~\ref{teo:meyers-serrin-haz-frontera}.
In either case there is a sequence
$(\mathbf u_k)_{k\in\mathbb N}\subseteq
\Gamma(\mathbf E)\cap W^{m,p}(M,\mathbf E)$ such that
$\mathbf u_k\to\mathbf u$ in $W^{m,p}(M,\mathbf E)$.
The estimate for smooth sections gives
\[
 \|\psi\mathbf u_k-\psi\mathbf u_\ell\|_{W^{m,p}(M,\mathbf E)}
 \leq C_{\psi,m}
       \|\mathbf u_k-\mathbf u_\ell\|_{W^{m,p}(M,\mathbf E)}.
\]
Thus $(\psi\mathbf u_k)_k$ is Cauchy and, by completeness of
$W^{m,p}(M,\mathbf E)$, has a limit
$\mathbf v\in W^{m,p}(M,\mathbf E)$.
The Sobolev norm controls the order-zero norm, so
$\psi\mathbf u_k\to\mathbf v$ also in $L^p(M,\mathbf E)$.
On the other hand, the pointwise product $\psi\mathbf u$ is
well-defined as a class of measurable sections and satisfies
\[
 \|\psi\mathbf u_k-\psi\mathbf u\|_{L^p(M,\mathbf E)}
 \leq\|\psi\|_{L^\infty(M)}
       \|\mathbf u_k-\mathbf u\|_{L^p(M,\mathbf E)}
 \longrightarrow0.
\]
Uniqueness of the limit in $L^p(M,\mathbf E)$ gives
$\mathbf v=\psi\mathbf u$ almost everywhere. In particular, the
pointwise product belongs to $W^{m,p}(M,\mathbf E)$ and the resulting
limit is independent of the approximating sequence.
Passing to the limit in the smooth estimate and using continuity
of the norms gives
\[
 \|\psi\mathbf u\|_{W^{m,p}(M,\mathbf E)}
 =\lim_{k\to\infty}\|\psi\mathbf u_k\|_{W^{m,p}(M,\mathbf E)}
 \leq C_{\psi,m}\|\mathbf u\|_{W^{m,p}(M,\mathbf E)}.
\]
Linearity follows from the pointwise product by the fixed function
$\psi$, and this bound proves continuity.
\end{proof} \begin{definition}[Sobolev spaces with zero trace]\label{def:W0mp-frontera}\index{Sobolev space@Sobolev space!with zero trace} Let $(M,\mathbf{g})$ be a Riemannian manifold with nonempty boundary and $\mathbf{E}\to M$ a smooth vector bundle equipped with a bundle metric $\mathbf{h}_{\mathbf{E}}$ and a smooth connection $\nabla^{\mathbf{E}}$. Let $m\in\mathbb N$ and $1\leq p<\infty$.

Define the Sobolev space with zero trace by \[
W^{m,p}_{0}(M,\mathbf{E})
:=
\overline{\Gamma_{c,\operatorname{Int}(M)}(\mathbf{E})}^{\,\|\cdot\|_{W^{m,p}(M,\mathbf{E})}}.
\]

That is, $W^{m,p}_{0}(M,\mathbf{E})$ is the closure, in the Sobolev norm, of the space of smooth sections with compact support contained in the interior of $M$. \end{definition} Before formulating the trace theorem, we must fix the connection to be used on the bundle restricted to the boundary. If $\mathbf{u}$ is a smooth section of $\mathbf{E}$, then $\mathbf{u}\restriction_{\partial M}$ is a smooth section of $\mathbf{E}\restriction_{\partial M}$. To measure the Sobolev regularity of this restriction, it does not suffice to regard $\mathbf{E}\restriction_{\partial M}$ simply as a vector bundle: we also need a connection on it. The natural connection is obtained by restricting the connection of $\mathbf{E}$ to directions tangent to the boundary.

\begin{lemma}[Induced connection on the restricted bundle]\label{lem:conexion-inducida-haz-restringido-frontera}\index{induced connection on the restricted bundle@induced connection on the restricted bundle} Let $(M,\mathbf{g})$ be a Riemannian manifold with nonempty boundary, and let $\mathbf{E}\to M$ be a smooth vector bundle equipped with a connection $\nabla^{\mathbf{E}}$. Then there is a connection on the restricted bundle $\mathbf{E}\restriction_{\partial M}\to\partial M$, denoted by $\nabla^{\mathbf{E}\restriction_{\partial M}}$, defined as follows.

If $\mathbf{X}\in\mathfrak{X}(\partial M)$ and $\mathbf{u}\in\Gamma(\mathbf{E}\restriction_{\partial M})$, take smooth extensions $\widetilde{\mathbf{X}}\in\mathfrak{X}(M)$ and $\widetilde{\mathbf{u}}\in\Gamma(\mathbf{E})$ such that $\widetilde{\mathbf{X}}\restriction_{\partial M}=\mathbf{X}$ and $\widetilde{\mathbf{u}}\restriction_{\partial M}=\mathbf{u}$, and define \[
\nabla^{\mathbf{E}\restriction_{\partial M}}_{\mathbf{X}}\mathbf{u}
:=
(\nabla^{\mathbf{E}}_{\widetilde{\mathbf{X}}}\widetilde{\mathbf{u}})\restriction_{\partial M}.
\] This definition is independent of the chosen extensions. Moreover, if $\nabla^{\mathbf{E}}$ is compatible with the bundle metric $\mathbf{h}_{\mathbf{E}}$, then $\nabla^{\mathbf{E}\restriction_{\partial M}}$ is compatible with the restricted bundle metric $\mathbf{h}_{\mathbf{E}}\restriction_{\mathbf{E}\restriction_{\partial M}}$. \end{lemma}

\begin{proof} We first prove independence of the chosen extensions. Let $\widetilde{\mathbf{X}}_1,\widetilde{\mathbf{X}}_2\in\mathfrak{X}(M)$ be two extensions of $\mathbf{X}$ and let $\widetilde{\mathbf{u}}_1,\widetilde{\mathbf{u}}_2\in\Gamma(\mathbf{E})$ be two extensions of $\mathbf{u}$. We want to prove that \[
(\nabla^{\mathbf{E}}_{\widetilde{\mathbf{X}}_1}\widetilde{\mathbf{u}}_1)\restriction_{\partial M}
=
(\nabla^{\mathbf{E}}_{\widetilde{\mathbf{X}}_2}\widetilde{\mathbf{u}}_2)\restriction_{\partial M}.
\] By linearity, \[
\nabla^{\mathbf{E}}_{\widetilde{\mathbf{X}}_1}\widetilde{\mathbf{u}}_1-
\nabla^{\mathbf{E}}_{\widetilde{\mathbf{X}}_2}\widetilde{\mathbf{u}}_2
=
\nabla^{\mathbf{E}}_{\widetilde{\mathbf{X}}_1-\widetilde{\mathbf{X}}_2}\widetilde{\mathbf{u}}_1
+
\nabla^{\mathbf{E}}_{\widetilde{\mathbf{X}}_2}(\widetilde{\mathbf{u}}_1-\widetilde{\mathbf{u}}_2).
\] The first term vanishes upon restriction to $\partial M$, since $\widetilde{\mathbf{X}}_1-\widetilde{\mathbf{X}}_2$ vanishes on $\partial M$ and the connection is $C^\infty(M)$-linear in its vector argument. For each $x\in\partial M$, the value of $\nabla^{\mathbf{E}}_{\widetilde{\mathbf{X}}_1-\widetilde{\mathbf{X}}_2}\widetilde{\mathbf{u}}_1$ at $x$ depends only on $(\widetilde{\mathbf{X}}_1-\widetilde{\mathbf{X}}_2)_x$, which is zero.

We now show that the second term also vanishes upon restriction to $\partial M$. Set $\mathbf{w}:=\widetilde{\mathbf{u}}_1-\widetilde{\mathbf{u}}_2$. Then $\mathbf{w}\restriction_{\partial M}=0$. Since the argument is local, take a local frame $(\mathbf{e}_1,\dots,\mathbf{e}_r)$ of $\mathbf{E}$ near a point of $\partial M$ and write $\mathbf{w}=\displaystyle\sum_{a=1}^{r}w^a\mathbf{e}_a$. Since $\mathbf{w}$ vanishes on $\partial M$ and the vectors $\mathbf{e}_1(x),\dots,\mathbf{e}_r(x)$ form a basis of $\mathbf{E}_x$ for each $x$, we have $w^a\restriction_{\partial M}=0$ for every $a$.

Now, using the Leibniz rule for the connection, \[
\nabla^{\mathbf{E}}_{\widetilde{\mathbf{X}}_2}\mathbf{w}
=
\displaystyle\sum_{a=1}^{r}\widetilde{\mathbf{X}}_2(w^a)\mathbf{e}_a
+
\displaystyle\sum_{a=1}^{r}w^a\nabla^{\mathbf{E}}_{\widetilde{\mathbf{X}}_2}\mathbf{e}_a.
\] Upon restriction to $\partial M$, the second summand vanishes because $w^a\restriction_{\partial M}=0$. For the first summand, use $\widetilde{\mathbf{X}}_2\restriction_{\partial M}=\mathbf{X}$ and the fact that $\mathbf{X}$ is tangent to $\partial M$. Since $w^a\restriction_{\partial M}=0$, we have $\mathbf{X}(w^a\restriction_{\partial M})=0$, and therefore $\widetilde{\mathbf{X}}_2(w^a)\restriction_{\partial M}=0$. Thus $(\nabla^{\mathbf{E}}_{\widetilde{\mathbf{X}}_2}\mathbf{w})\restriction_{\partial M}=0$. This proves that $\nabla^{\mathbf{E}\restriction_{\partial M}}_{\mathbf{X}}\mathbf{u}$ is well defined.

We now prove that $\nabla^{\mathbf{E}\restriction_{\partial M}}$ is a connection. Let $f\in C^\infty(\partial M)$, $\mathbf{X},\mathbf{Y}\in\mathfrak{X}(\partial M)$, and $\mathbf{u},\mathbf{v}\in\Gamma(\mathbf{E}\restriction_{\partial M})$. Take smooth extensions of these objects to a neighborhood of $\partial M$. Since $\nabla^{\mathbf{E}}$ is $C^\infty(M)$-linear in its vector argument, we obtain $\nabla^{\mathbf{E}\restriction_{\partial M}}_{f\mathbf{X}+\mathbf{Y}}\mathbf{u}=f\nabla^{\mathbf{E}\restriction_{\partial M}}_{\mathbf{X}}\mathbf{u}+\nabla^{\mathbf{E}\restriction_{\partial M}}_{\mathbf{Y}}\mathbf{u}$. Likewise, using the Leibniz rule in the section argument, \[
\nabla^{\mathbf{E}\restriction_{\partial M}}_{\mathbf{X}}(f\mathbf{u})
=
\mathbf{X}(f)\mathbf{u}+f\nabla^{\mathbf{E}\restriction_{\partial M}}_{\mathbf{X}}\mathbf{u}.
\] Therefore, $\nabla^{\mathbf{E}\restriction_{\partial M}}$ is a connection on $\mathbf{E}\restriction_{\partial M}$.

Finally, suppose that $\nabla^{\mathbf{E}}$ is compatible with $\mathbf{h}_{\mathbf{E}}$. Let $\mathbf{X}\in\mathfrak{X}(\partial M)$ and $\mathbf{u},\mathbf{v}\in\Gamma(\mathbf{E}\restriction_{\partial M})$. Take extensions $\widetilde{\mathbf{X}},\widetilde{\mathbf{u}},\widetilde{\mathbf{v}}$ as before. By compatibility of $\nabla^{\mathbf{E}}$ with $\mathbf{h}_{\mathbf{E}}$, \[
\widetilde{\mathbf{X}}\big(\mathbf{h}_{\mathbf{E}}(\widetilde{\mathbf{u}},\widetilde{\mathbf{v}})\big)
=
\mathbf{h}_{\mathbf{E}}(\nabla^{\mathbf{E}}_{\widetilde{\mathbf{X}}}\widetilde{\mathbf{u}},\widetilde{\mathbf{v}})
+
\mathbf{h}_{\mathbf{E}}(\widetilde{\mathbf{u}},\nabla^{\mathbf{E}}_{\widetilde{\mathbf{X}}}\widetilde{\mathbf{v}}).
\] Upon restriction to $\partial M$, the left-hand side becomes $\mathbf{X}\big(\mathbf{h}_{\mathbf{E}}\restriction_{\mathbf{E}\restriction_{\partial M}}(\mathbf{u},\mathbf{v})\big)$, while the right-hand side becomes \[
\mathbf{h}_{\mathbf{E}}\restriction_{\mathbf{E}\restriction_{\partial M}}(\nabla^{\mathbf{E}\restriction_{\partial M}}_{\mathbf{X}}\mathbf{u},\mathbf{v})
+
\mathbf{h}_{\mathbf{E}}\restriction_{\mathbf{E}\restriction_{\partial M}}(\mathbf{u},\nabla^{\mathbf{E}\restriction_{\partial M}}_{\mathbf{X}}\mathbf{v}).
\] Therefore, $\nabla^{\mathbf{E}\restriction_{\partial M}}$ is compatible with the restricted bundle metric. \end{proof}

Henceforth, when we write $W^{m,p}(\partial M,\mathbf{E}\restriction_{\partial M})$, we understand that this space is constructed on the Riemannian manifold $\partial M$ using the restricted bundle $\mathbf{E}\restriction_{\partial M}$, the restricted bundle metric $\mathbf{h}_{\mathbf{E}}\restriction_{\mathbf{E}\restriction_{\partial M}}$, and the induced connection $\nabla^{\mathbf{E}\restriction_{\partial M}}$.

\subsection{The trace operator}

In a boundary chart, restriction of a section to $\partial M$ is expressed through the trace of its components on the flat boundary. We use Theorem~\ref{teo:traza-semiespacio}, Lemma~\ref{lem:extension-cero-media-bola-semiespacio}, and Corollary~\ref{cor:caracterizacion-local-media-bola}, proved in the section on Euclidean traces in Chapter~\ref{cap:sobolev-euclidiano}. In particular, components localized in a half-ball are extended by zero within $\mathbb R^n_+$ before taking the trace. We retain the convention fixed there for writing $\gamma_0D^\beta v$ on the flat face of a half-ball.

Before stating the global result, recall that, when we write \(W^{k,p}(\partial M,\mathbf{E}\restriction_{\partial M})\), we use the induced connection on the restricted bundle \(\mathbf{E}\restriction_{\partial M}\to\partial M\), constructed in Lemma~\ref{lem:conexion-inducida-haz-restringido-frontera}. This matters because the trace of a section in \(\mathbf{E}\) will take values in a Sobolev space on \(\mathbf{E}\restriction_{\partial M}\).

\begin{theorem}[Integer-order trace theorem for vector bundles]\label{teo:traza-entera-haces}\index{trace operator!of integer order for vector bundles}\glsadd{operador-traza} Let \((M,\mathbf{g})\) be a compact Riemannian manifold with nonempty boundary, and let \(\mathbf{E}\to M\) be a smooth vector bundle of finite rank, equipped with a bundle metric \(\mathbf{h}_{\mathbf{E}}\) and a compatible connection \(\nabla^{\mathbf{E}}\). Let \(m\geq1\) and \(1\leq p<\infty\). Then the restriction operator \[
\gamma_0\colon \Gamma(\mathbf{E})\longrightarrow\Gamma(\mathbf{E}\restriction_{\partial M}),
\qquad
\gamma_0(\mathbf{u})=\mathbf{u}\restriction_{\partial M},
\] admits a unique continuous linear extension \[
\boldsymbol{\operatorname{Tr}}\colon W^{m,p}(M,\mathbf{E})\longrightarrow W^{m-1,p}(\partial M,\mathbf{E}\restriction_{\partial M}).
\] In particular, there is a constant \(C>0\) such that \[
\|\boldsymbol{\operatorname{Tr}}(\mathbf{u})\|_{W^{m-1,p}(\partial M,\mathbf{E}\restriction_{\partial M})}
\leq
C\|\mathbf{u}\|_{W^{m,p}(M,\mathbf{E})}
\] for every \(\mathbf{u}\in W^{m,p}(M,\mathbf{E})\). \end{theorem}

\begin{proof} We first prove the estimate for smooth sections. Since \(\partial M\) is compact, take a finite family of regular coordinate half-balls \((U_\alpha,\phi_\alpha)_{\alpha=1}^{N}\) covering \(\partial M\). Write \(B_\alpha^+:=\phi_\alpha(U_\alpha\cap\operatorname{Int}(M))\) and \(B_\alpha':=\phi_\alpha(U_\alpha\cap\partial M)\). Refining if necessary, we may assume that $\overline{U_\alpha}$ lies in the domain of a smooth local frame \((\mathbf{e}_{1,\alpha},\dots,\mathbf{e}_{r,\alpha})\) trivializing \(\mathbf{E}\). Take a smooth partition of unity \((\psi_\alpha)_{\alpha=1}^{N}\), defined on a neighborhood of \(\partial M\), subordinate to the charts \(U_\alpha\) and such that \(\operatorname{supp}\psi_\alpha\Subset U_\alpha\).

Let \(\mathbf{u}\in\Gamma(\mathbf{E})\). On \(U_\alpha\) write \(\mathbf{u}=\displaystyle\sum_{a=1}^{r}u_\alpha^a\mathbf{e}_{a,\alpha}\). Then \(\psi_\alpha \mathbf{u}=\displaystyle\sum_{a=1}^{r}\psi_\alpha u_\alpha^a\mathbf{e}_{a,\alpha}\). For each component, define \(f_{\alpha}^{a}:=(\psi_\alpha u_\alpha^a)\circ\phi_\alpha^{-1}\) on \(B_\alpha^+\). Since \(\operatorname{supp}\psi_\alpha\Subset U_\alpha\), the function \(f_\alpha^a\) vanishes on a neighborhood of the spherical part of \(B_\alpha^+\). Let \(\widetilde f_\alpha^a\) be its extension by zero to \(\mathbb R^n_+\). By Lemma~\ref{lem:extension-cero-media-bola-semiespacio}, \(\widetilde f_\alpha^a\in W^{m,p}(\mathbb R^n_+\)), and moreover \(\|\widetilde f_\alpha^a\|_{W^{m,p}(\mathbb R^n_+)}=\|f_\alpha^a\|_{W^{m,p}(B_\alpha^+)}\). Applying Theorem~\ref{teo:traza-semiespacio} to \(\widetilde f_\alpha^a\), and using continuity of restriction from \(W^{m-1,p}(\mathbb R^{n-1})\) to \(W^{m-1,p}(B_\alpha')\), gives \[
\left\|
(\gamma_0\widetilde f_\alpha^a)\restriction_{B_\alpha'}
\right\|_{W^{m-1,p}(B_\alpha')}
\leq
C_\alpha
\left\|
f_\alpha^a
\right\|_{W^{m,p}(B_\alpha^+)}.
\] For smooth \(\mathbf{u}\), \((\gamma_0\widetilde f_\alpha^a)\restriction_{B_\alpha'}\) agrees with \((f_\alpha^a)\restriction_{B_\alpha'}\), since \(\widetilde f_\alpha^a=f_\alpha^a\) near \(B_\alpha'\) within the half-space. Summing over \(a\) and using local equivalence of norms in the regular coordinate half-ball gives \[
\displaystyle\sum_{a=1}^{r}
\left\|
(\gamma_0\widetilde f_\alpha^a)\restriction_{B_\alpha'}
\right\|_{W^{m-1,p}(B_\alpha')}
\leq
C_\alpha'
\|\psi_\alpha \mathbf{u}\|_{W^{m,p}(U_\alpha,\mathbf{E}\restriction_{U_\alpha})}.
\] On the other hand, applying local equivalence of norms on the manifold without boundary \(\partial M\) to the restricted bundle \(\mathbf{E}\restriction_{\partial M}\to\partial M\) gives \[
\|(\psi_\alpha \mathbf{u})\restriction_{\partial M}\|_{W^{m-1,p}(U_\alpha\cap\partial M,\mathbf{E}\restriction_{U_\alpha\cap\partial M})}
\leq
C_\alpha''
\displaystyle\sum_{a=1}^{r}
\left\|
(\gamma_0\widetilde f_\alpha^a)\restriction_{B_\alpha'}
\right\|_{W^{m-1,p}(B_\alpha')}.
\] Therefore, \[
\|(\psi_\alpha \mathbf{u})\restriction_{\partial M}\|_{W^{m-1,p}(U_\alpha\cap\partial M,\mathbf{E}\restriction_{U_\alpha\cap\partial M})}
\leq
C_\alpha'''
\|\psi_\alpha \mathbf{u}\|_{W^{m,p}(U_\alpha,\mathbf{E}\restriction_{U_\alpha})}.
\] Now using Lemma~\ref{lem:multiplicacion-funcion-suave-sobolev}, applied to the fixed function \(\psi_\alpha\), there is a constant \(C_\alpha''''>0\) such that \(\|\psi_\alpha \mathbf{u}\|_{W^{m,p}(M,\mathbf{E})}\leq C_\alpha''''\|\mathbf{u}\|_{W^{m,p}(M,\mathbf{E})}\). Since restricting the global norm to \(U_\alpha\) gives a value no greater than the global norm, we conclude that \[
\|(\psi_\alpha \mathbf{u})\restriction_{\partial M}\|_{W^{m-1,p}(U_\alpha\cap\partial M,\mathbf{E}\restriction_{U_\alpha\cap\partial M})}
\leq
C_\alpha'''''
\|\mathbf{u}\|_{W^{m,p}(M,\mathbf{E})}.
\]

Since there are finitely many charts and \(\displaystyle\sum_{\alpha=1}^{N}\psi_\alpha=1\) on a neighborhood of \(\partial M\), we have \(\mathbf{u}\restriction_{\partial M}=\displaystyle\sum_{\alpha=1}^{N}(\psi_\alpha \mathbf{u})\restriction_{\partial M}\). Using the triangle inequality and local equivalence of norms on \(\mathbf{E}\restriction_{\partial M}\) gives a constant \(C>0\) such that \[
\|\mathbf{u}\restriction_{\partial M}\|_{W^{m-1,p}(\partial M,\mathbf{E}\restriction_{\partial M})}
\leq
C\|\mathbf{u}\|_{W^{m,p}(M,\mathbf{E})}
\] for every \(\mathbf{u}\in\Gamma(\mathbf{E})\).

The preceding estimate shows that restriction is a continuous linear operator from $\Gamma(\mathbf{E})\cap W^{m,p}(M,\mathbf{E})$ to $W^{m-1,p}(\partial M,\mathbf{E}\restriction_{\partial M})$. Theorem~\ref{teo:meyers-serrin-haz-frontera} states that its domain is dense in $W^{m,p}(M,\mathbf{E})$, while the target space is Banach. Theorem~\ref{teo:extension-operadores-subespacio-denso} therefore gives a unique continuous linear extension \[
\boldsymbol{\operatorname{Tr}}\colon W^{m,p}(M,\mathbf{E})\longrightarrow W^{m-1,p}(\partial M,\mathbf{E}\restriction_{\partial M})
\] and preserves the estimate \[
\|\boldsymbol{\operatorname{Tr}}(\mathbf{u})\|_{W^{m-1,p}(\partial M,\mathbf{E}\restriction_{\partial M})}
\leq
C\|\mathbf{u}\|_{W^{m,p}(M,\mathbf{E})}.
\] \end{proof} \begin{figura}[htbp] \includegraphics[width=\linewidth]{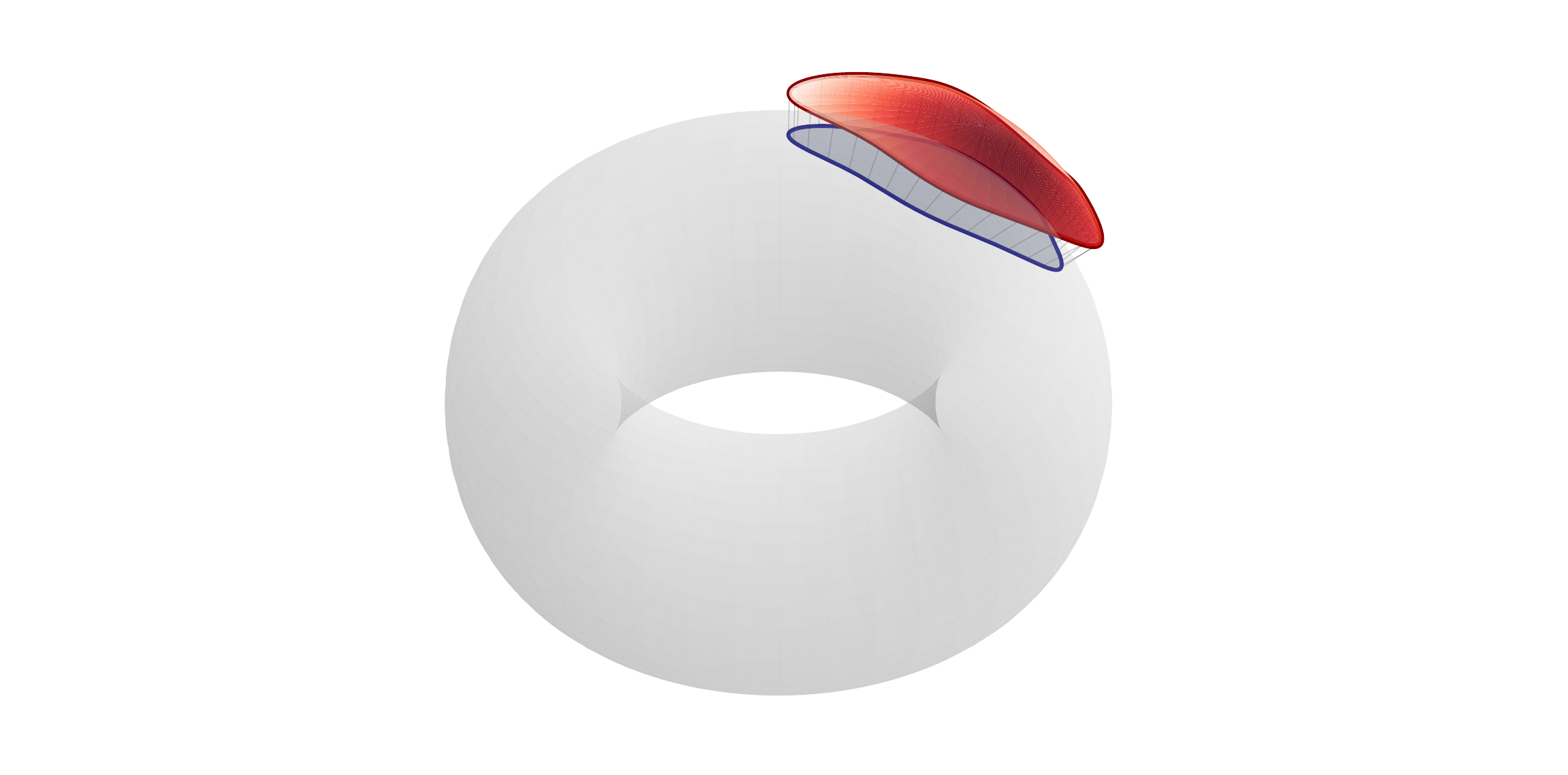} \caption{The trace theorem generalizes the idea of restricting a smooth function to the boundary of a manifold. Here \(\Omega\subseteq T^2\) is an open subset of the torus whose closure \(\overline{\Omega}\) is a manifold with boundary. The blue region represents \(\Omega\), while its boundary \(\partial\Omega\) is shown in dark blue. The graph of the function \(u\colon \overline{\Omega}\longrightarrow\mathbb{R}\) is shown in red and its trace on \(\partial\Omega\) in dark red.} \label{fig:traza-en-el-toro} \end{figura} \subsection{Intrinsic Sobolev--Slobodeckij spaces and the exact trace}

The integer-order trace theorem gives an estimate in $W^{m-1,p}(\partial M,\mathbf{E}\restriction_{\partial M})$, but this regularity does not fully describe the image. Restriction to a hypersurface loses exactly $\frac{1}{p}$ derivatives. To formulate this assertion without choosing charts or frames, we compare values of a section in different fibers by parallel transport for the induced connection on the boundary.

The construction is presented first on a closed manifold of positive dimension, since $\partial M$ is such a manifold. Let $N$ be a compact Riemannian manifold without boundary, and let $\mathbf{F}\longrightarrow N$ be a smooth vector bundle with a bundle metric and a compatible connection. If $N$ is disconnected, the following definitions apply to the restrictions of these data to each component; compactness ensures that there are only finitely many. For $s>0$ and $1\leq p<\infty$, if $N=N_1\sqcup\cdots\sqcup N_q$, the global norm is defined by \[
\|\mathbf{v}\|_{W^{s,p}(N,\mathbf{F})}
:=
\left(\sum_{\ell=1}^q
\|\mathbf{v}\restriction_{N_\ell}\|_{W^{s,p}(N_\ell,\mathbf{F}\restriction_{N_\ell})}^p\right)^{\frac{1}{p}}.
\] In dimension zero, set $W^{s,p}(N,\mathbf{F}):=L^p(N,\mathbf{F})$ for every $s>0$; this is then a finite sum of fibers, with no difference seminorm.

\begin{definition}[Intrinsic Gagliardo seminorm for sections] \label{def:seminorma-gagliardo-intrinseca-secciones} \index{Gagliardo seminorm!for sections} \index{Slobodeckij space!for sections} Let $(N,\mathbf{h})$ be a compact connected Riemannian manifold without boundary of dimension $d\geq1$, and let $\mathbf{F}\longrightarrow N$ be a smooth vector bundle equipped with a bundle metric $\mathbf{h}_{\mathbf{F}}$ and a compatible connection $\nabla^{\mathbf{F}}$. Set \[
\rho_N:=\min\{1,\operatorname{inj}(N)\}.
\] The compact manifold $N$ is complete, and Proposition~\ref{prop: continuidad radio inyectividad conexa y completa} shows that the pointwise injectivity radius is continuous and positive. Its minimum on $N$ is therefore strictly positive. For $0<\sigma<1$, $1\leq p<\infty$, and $\mathbf{v}\in L^p(N,\mathbf{F})$ define \begin{equation}
\label{eq:seminorma-gagliardo-intrinseca-secciones}
[\mathbf{v}]_{W^{\sigma,p}(N,\mathbf{F})}
:=
\left(
\int_N
\int_{\{y\in N\mid\,0<d_{\mathbf{h}}(x,y)<\rho_N\}}
\frac{|\mathbf{v}(x)-\mathbf{v}(y)|_{\mathbf{h}_{\mathbf{F}};x,y}^{p}}
{d_{\mathbf{h}}(x,y)^{d+\sigma p}}\,d\lambda_{\mathbf{h}}(y)\,d\lambda_{\mathbf{h}}(x)
\right)^{\frac{1}{p}}.
\end{equation} The fiberwise separation in the numerator is that of Definition~\ref{def:separacion-fibrada-transporte-paralelo}; thus \[
|\mathbf{v}(x)-\mathbf{v}(y)|_{\mathbf{h}_{\mathbf{F}};x,y}
=
\left|\mathbf{v}(x)-P^{\mathbf{F}}_{\gamma_{xy},1\to0}\mathbf{v}(y)\right|_{\mathbf{h}_{\mathbf{F}}(x)},
\] where $\gamma_{xy}$ is the unique short minimizing geodesic from $x$ to $y$. Transport depends smoothly on $(x,y)$ in the region $0<d_{\mathbf{h}}(x,y)<\rho_N$, so the integrand is measurable. Changing a representative on a null set $Z\subseteq N$ changes the integrand only on $(Z\times N)\cup(N\times Z)$, which has zero product measure by Theorem~\ref{teo:tonelli}.

Define \[
W^{\sigma,p}(N,\mathbf{F})
:=
\left\{\mathbf{v}\in L^p(N,\mathbf{F})\middle|
[\mathbf{v}]_{W^{\sigma,p}(N,\mathbf{F})}<\infty\right\}
\] and equip it with the norm \[
\|\mathbf{v}\|_{W^{\sigma,p}(N,\mathbf{F})}
:=
\|\mathbf{v}\|_{L^p(N,\mathbf{F})}+[\mathbf{v}]_{W^{\sigma,p}(N,\mathbf{F})}.
\]

If $s=k+\sigma$, with $k\in\mathbb N_0$ and $0<\sigma<1$, set \[
\mathbf{F}_k:=T^{(0,k)}(TN)\otimes \mathbf{F}
\] with its product metric and connection, and define \begin{equation}
\label{eq:definicion-slobodeckij-intrinseco-orden-superior}
W^{s,p}(N,\mathbf{F})
:=
\left\{\mathbf{v}\in W^{k,p}(N,\mathbf{F})\middle|
[\nabla_w^k \mathbf{v}]_{W^{\sigma,p}(N,\mathbf{F}_k)}<\infty\right\},
\end{equation} with \begin{equation}
\label{eq:norma-slobodeckij-intrinseca-orden-superior}
\|\mathbf{v}\|_{W^{s,p}(N,\mathbf{F})}
:=
\|\mathbf{v}\|_{W^{k,p}(N,\mathbf{F})}
+[\nabla_w^k \mathbf{v}]_{W^{\sigma,p}(N,\mathbf{F}_k)}.
\end{equation} For $s\in\mathbb N_0$, retain the preceding intrinsic definition of $W^{s,p}(N,\mathbf{F})$ using weak covariant derivatives. \end{definition}

The seminorm of order $0<\sigma<1$ is the version, for a fixed smooth connection, of the covariant fractional spaces considered in \cite[Definition~6.4]{VanSchaftingenWinter2025}; higher orders are obtained by applying it to $\nabla_w^k \mathbf{v}$ with the product connection. Its relationship with Euclidean theory is established by the following localized characterization.

\begin{lemma}[Independence of the comparison radius] \label{lem:independencia-radio-gagliardo-intrinseco} Let $(N,\mathbf{h})$ be a compact connected Riemannian manifold without boundary of dimension $d\geq1$, and let $\mathbf{F}\to N$ be a smooth vector bundle with a bundle metric and a compatible connection. Let $0<\sigma<1$, $1\leq p<\infty$, and $0<r_1\leq r_2\leq\operatorname{inj}(N)$. Denote by $[\mathbf{v}]_{\sigma,p;r_i}$ the expression in \eqref{eq:seminorma-gagliardo-intrinseca-secciones} with $r_i$ in place of $\rho_N$. There are constants $c,C>0$, independent of $v$, such that \[
c\bigl(\|\mathbf{v}\|_{L^p(N,\mathbf{F})}+[\mathbf{v}]_{\sigma,p;r_1}\bigr)
\leq
\|\mathbf{v}\|_{L^p(N,\mathbf{F})}+[\mathbf{v}]_{\sigma,p;r_2}
\leq
C\bigl(\|\mathbf{v}\|_{L^p(N,\mathbf{F})}+[\mathbf{v}]_{\sigma,p;r_1}\bigr).
\] \end{lemma}

\begin{proof} The first inequality follows from inclusion of the regions of integration. In the region $r_1\leq d_{\mathbf{h}}(x,y)<r_2$, the isometry property of parallel transport and the inequality $|a-b|^p\leq2^{p-1}(|a|^p+|b|^p)$ give \[
\frac{|\mathbf{v}(x)-\mathbf{v}(y)|_{\mathbf{h}_{\mathbf{F}};x,y}^p}{d_{\mathbf{h}}(x,y)^{d+\sigma p}}
\leq
2^{p-1}r_1^{-d-\sigma p}
\bigl(|\mathbf{v}(x)|_{\mathbf{h}_{\mathbf{F}}}^p+|\mathbf{v}(y)|_{\mathbf{h}_{\mathbf{F}}}^p\bigr).
\] Integrating over $N\times N$ and using $\lambda_{\mathbf{h}}(N)<\infty$ bounds the contribution of this region by $C\|\mathbf{v}\|_{L^p(N,\mathbf{F})}^p$. This proves the second inequality. \end{proof}

The following Euclidean fact will control the lower-order terms that arise when a covariant derivative is expressed in a frame.

\begin{lemma}[Fractional control by one integer-order derivative] \label{lem:W1p-controla-slobodeckij-global} Let $0<\sigma<1$ and $1\leq p<\infty$. There is $C>0$ such that \[
\|f\|_{W^{\sigma,p}(\mathbb R^d)}
\leq C\|f\|_{W^{1,p}(\mathbb R^d)}
\] for every $f\in W^{1,p}(\mathbb R^d)$. Moreover, if $a\in C_c^\infty(\mathbb R^d)$, then \[
\|af\|_{W^{\sigma,p}(\mathbb R^d)}
\leq C_a\|f\|_{W^{\sigma,p}(\mathbb R^d)}.
\] \end{lemma}

\begin{proof} By Proposition~\ref{prop:gagliardo-mediante-incrementos} and Lemma~\ref{lem:incrementos-W1p}, for $\|z\|<1$ we have \[
\|\tau_zf-f\|_{L^p(\mathbb R^d)}^p
\leq \|z\|^p\|Df\|_{L^p(\mathbb R^d,\mathbb R^d)}^p,
\] while for $\|z\|\geq1$ we use $\|\tau_zf-f\|_{L^p(\mathbb R^d)}^p\leq2^p\|f\|_{L^p(\mathbb R^d)}^p$. The resulting radial integrals are \[
\int_0^1r^{p(1-\sigma)-1}\,dr<\infty,
\qquad
\int_1^\infty r^{-\sigma p-1}\,dr<\infty.
\] This proves the first assertion.

For the second, write \[
a(x)f(x)-a(y)f(y)
=a(x)(f(x)-f(y))+(a(x)-a(y))f(y).
\] The first term is controlled by $\|a\|_{L^\infty(\mathbb R^d)}[f]_{W^{\sigma,p}(\mathbb R^d)}$. For the second, use \[
|a(x)-a(y)|
\leq C_a\min\{1,\|x-y\|\};
\] the same separation into $\|x-y\|<1$ and $\|x-y\|\geq1$ gives a bound by $C_a\|f\|_{L^p(\mathbb R^d)}$. Together with the estimate $L^p(\mathbb R^d)$, this yields the result. \end{proof}

Such coordinate balls and frames can be chosen around every point:
apply Lemma~\ref{lem:vecindades-normales-convexas-locales} inside a
trivializing neighborhood and then shrink the radius and the inner ball.

\begin{lemma}[Local comparison of transport and components] \label{lem:comparacion-local-slobodeckij-transporte-componentes} Let $(N,\mathbf{h})$ be a compact connected Riemannian manifold without boundary, and let $\mathbf{F}\to N$ be a smooth vector bundle with a bundle metric and a compatible connection. Let $(U,\phi)$ be a regular coordinate ball in $N$. Suppose that there is a geodesically convex normal ball $U'$ such that $\overline U\subseteq U'$, $\operatorname{diam}_{\mathbf{h}}(U')<\rho_N$, and both $\phi$ and the smooth frame $(\mathbf{e}_1,\dots,\mathbf{e}_r)$ of $\mathbf{F}$ extend to $U'$. Write $\displaystyle \mathbf{v}=\displaystyle\sum_{a=1}^rv^a\mathbf{e}_a$ on $U$. If $0<\sigma<1$ and $1\leq p<\infty$, there are constants $c,C>0$, independent of $\mathbf{v}$, such that \begin{equation}
\label{eq:comparacion-local-slobodeckij-transporte-componentes}
\begin{aligned}
c\sum_{a=1}^r
\|v^a\circ\phi^{-1}\|_{W^{\sigma,p}(\phi(U))}
&\leq
\|\mathbf{v}\|_{L^p(U,\mathbf{F})}+[\mathbf{v}]_{W^{\sigma,p}(U,\mathbf{F})}\\
&\leq
C\sum_{a=1}^r
\|v^a\circ\phi^{-1}\|_{W^{\sigma,p}(\phi(U))},
\end{aligned}
\end{equation} where the middle seminorm integrates only over pairs of points in $U$. \end{lemma}

\begin{proof} Set $\Omega:=\phi(U)$. In the chart and frame, transport from $y$ to $x$ is represented by a matrix $G(x,y)$. If $c_{xy}:=\phi\circ\gamma_{xy}$, the local parallel transport equation, Proposition~\ref{prop: ecuacion local transporte paralelo operador}, is a linear equation whose coefficient matrices are the connection matrices evaluated along $\dot c_{xy}$. Those matrices are bounded on the compact set under consideration, and the length of $c_{xy}$ is comparable to $d_{\mathbf{h}}(x,y)$. The integral form of the equation and Grönwall's inequality in Lemma~\ref{lema: gronwall} therefore give \begin{equation}
\label{eq:transporte-local-cercano-identidad-slobodeckij}
\|G(x,y)-I\|
\leq C_0d_{\mathbf{h}}(x,y).
\end{equation} The same estimate holds for $G(x,y)^{-1}=G(y,x)$.

Uniform equivalence of the bundle norm and the component norm on $\overline U$ allows us to deduce from \eqref{eq:transporte-local-cercano-identidad-slobodeckij} the two inequalities \begin{align}
|\mathbf{v}(x)-\mathbf{v}(y)|_{\mathbf{h}_{\mathbf{F}};x,y}
&\leq C_1\bigl(|v^\flat(x)-v^\flat(y)|
+d_{\mathbf{h}}(x,y)|v^\flat(y)|\bigr),
\label{eq:transporte-a-componentes-slobodeckij}\\
|v^\flat(x)-v^\flat(y)|
&\leq C_1\bigl(|\mathbf{v}(x)-\mathbf{v}(y)|_{\mathbf{h}_{\mathbf{F}};x,y}
+d_{\mathbf{h}}(x,y)|v^\flat(y)|\bigr),
\label{eq:componentes-a-transporte-slobodeckij}
\end{align} where $v^\flat=(v^1,\dots,v^r)$. On a coordinate ball with compact closure there are constants $c_2,C_2>0$ such that \[
c_2\|\phi(x)-\phi(y)\|
\leq d_{\mathbf{h}}(x,y)
\leq C_2\|\phi(x)-\phi(y)\|,
\qquad
c_2\,dz\leq d\lambda_{\mathbf{h}}\leq C_2\,dz.
\] Divide the $p$th powers in \eqref{eq:transporte-a-componentes-slobodeckij} and \eqref{eq:componentes-a-transporte-slobodeckij} by $d_{\mathbf{h}}(x,y)^{d+\sigma p}$ and integrate. The error term is controlled by $\|\mathbf{v}\|_{L^p(U,\mathbf{F})}^p$, since, uniformly in $y$, \[
\int_{\phi(U)}
\frac{\|z-w\|^p}{\|z-w\|^{d+\sigma p}}\,dz
\leq
C\int_0^{C_2\operatorname{diam}(U)}r^{p(1-\sigma)-1}\,dr<\infty.
\] Comparison of the $L^p(N,\mathbf F)$ norms completes both inequalities in \eqref{eq:comparacion-local-slobodeckij-transporte-componentes}. \end{proof}

\begin{lemma}[Multiplication by smooth functions] \label{lem:multiplicacion-suave-slobodeckij-intrinseco} Let $(N,\mathbf{h})$ be a compact connected Riemannian manifold without boundary, and let $\mathbf{F}\to N$ be a smooth vector bundle with a bundle metric and a compatible connection. Let $s>0$, $s\notin\mathbb N$, $1\leq p<\infty$, and $\chi\in C^\infty(N)$. There is a constant $C_{\chi,s,p}>0$ such that \[
\|\chi\mathbf{v}\|_{W^{s,p}(N,\mathbf{F})}
\leq C_{\chi,s,p}\|\mathbf{v}\|_{W^{s,p}(N,\mathbf{F})}.
\] \end{lemma}

\begin{proof} First suppose that $0<s<1$. Compatibility of parallel transport with scalar multiplication gives \[
\begin{aligned}
|\chi(x)\mathbf{v}(x)-\chi(y)\mathbf{v}(y)|_{\mathbf{h}_{\mathbf{F}};x,y}
\leq{}&|\chi(x)|\,|\mathbf{v}(x)-\mathbf{v}(y)|_{\mathbf{h}_{\mathbf{F}};x,y}\\
&+|\chi(x)-\chi(y)|\,|\mathbf{v}(y)|_{\mathbf{h}_{\mathbf{F}}}.
\end{aligned}
\] Since $N$ is compact, $\chi$ is bounded and Lipschitz. The second term is integrated using $\displaystyle \int_0^{\rho_N}r^{p(1-s)-1}\,dr<\infty$, exactly as in the preceding proof.

Now let $s=k+\sigma$, with $k\geq1$. The iterated Leibniz rule expresses $\nabla_w^k(\chi\mathbf{v})$ as a finite sum of contractions of $\nabla^{k-j}\chi$ with $\nabla_w^j\mathbf{v}$. The term $j=k$ is controlled by the case of order $\sigma$ applied to $\nabla_w^k\mathbf{v}$. If $j<k$, then $\nabla_w^j\mathbf{v}\in W^{1,p}(N,\mathbf{F}_j)$. Let us justify this assertion and, at the same time, the Leibniz rule just used. By Theorem~\ref{meyers-serrin-haz}, there is a sequence of smooth sections $(\mathbf{v}_\nu)$ converging to $\mathbf{v}$ in $W^{k,p}(N,\mathbf{F})$. For $j<k$, the identities \[
\nabla^j\mathbf{v}_\nu\longrightarrow\nabla_w^j\mathbf{v},
\qquad
\nabla(\nabla^j\mathbf{v}_\nu)=\nabla^{j+1}\mathbf{v}_\nu
\longrightarrow\nabla_w^{j+1}\mathbf{v}
\] hold in the corresponding spaces $L^p(N,\mathbf F_j\oplus\mathbf F_{j+1})$. Lemma~\ref{derivada debil es operador cerrado}, applied to the bundle $\mathbf{F}_j$, shows that $\nabla_w^j\mathbf{v}$ belongs to $W^{1,p}(N,\mathbf{F}_j)$ and its first weak derivative is $\nabla_w^{j+1}\mathbf{v}$. The Leibniz formula for $\mathbf{v}_\nu$ passes to the limit in $L^p(N,\mathbf F_k)$ because its coefficients are the smooth bounded derivatives of $\chi$; this establishes the same formula for $\mathbf{v}$.

Now localize in a finite cover and apply a fixed cutoff equal to one on the support of each localization. The cut-off components, extended by zero, belong to $W^{1,p}(\mathbb R^d)$; Lemma~\ref{lem:W1p-controla-slobodeckij-global} then applies on $\mathbb R^d$. The cross term from the extension is controlled by the positive distance of the support from the chart boundary. Lemma~\ref{lem:comparacion-local-slobodeckij-transporte-componentes} allows us to return to the intrinsic seminorm. In each frame, contractions with $\nabla^{k-j}\chi$ are finite sums of multiplications by smooth coefficients; the second part of the same Euclidean lemma controls them. The derivatives of $\chi$ are bounded and the family of charts is finite, so every term is controlled by $C\|\mathbf{v}\|_{W^{s,p}(N,\mathbf{F})}$. The integer-order part of the norm is controlled directly by the same Leibniz rule, the uniform bounds on the derivatives of $\chi$, and the definition of $W^{k,p}(N,\mathbf{F})$. This completes the estimate of the integer-order norm and the fractional seminorm appearing in \eqref{eq:norma-slobodeckij-intrinseca-orden-superior}. \end{proof}

\begin{proposition}[Localized characterization of the intrinsic norm] \label{prop:caracterizacion-local-slobodeckij-intrinseco} Let $(N,\mathbf{h})$ be a compact connected Riemannian manifold without boundary of dimension $d\geq1$, and let $\mathbf{F}\to N$ be a smooth vector bundle of rank $r$ with a bundle metric and a compatible connection. Let $s>0$, $s\notin\mathbb N$, and $1\leq p<\infty$. Choose a finite cover $(U_i,\phi_i)_{i=1}^N$ of $N$ by regular coordinate balls as in Lemma~\ref{lem:comparacion-local-slobodeckij-transporte-componentes}, a partition of unity $(\psi_i)_{i=1}^N$ with $\operatorname{supp}\psi_i\Subset U_i$, and frames $(\mathbf{e}_{i,1},\dots,\mathbf{e}_{i,r})$ defined on neighborhoods of $\overline{U_i}$. Write \[
\psi_i\mathbf{v}=\sum_{a=1}^rv_i^a\mathbf{e}_{i,a}
\qquad\text{in }U_i
\] and denote by $\widetilde v_i^a$ the extension of $v_i^a\circ\phi_i^{-1}$ by zero to $\mathbb R^d$. There are constants $c,C>0$, independent of $\mathbf{v}$, such that \begin{equation}
\label{eq:caracterizacion-local-slobodeckij-intrinseco}
c\|\mathbf{v}\|_{W^{s,p}(N,\mathbf{F})}
\leq
\sum_{i=1}^N\sum_{a=1}^r
\|\widetilde v_i^a\|_{W^{s,p}(\mathbb R^d)}
\leq
C\|\mathbf{v}\|_{W^{s,p}(N,\mathbf{F})}.
\end{equation} In particular, the intrinsically defined space is independent of the charts, frames, and partition used to compute the localized norm. \end{proposition}

\begin{proof} First consider $0<s<1$. Lemma~\ref{lem:multiplicacion-suave-slobodeckij-intrinseco} controls $\psi_i\mathbf{v}$ by $\mathbf{v}$, and Lemma~\ref{lem:comparacion-local-slobodeckij-transporte-componentes} controls its components inside $\phi_i(U_i)$. Since $K_i:=\phi_i(\operatorname{supp}\psi_i)\Subset\phi_i(U_i)$, extension by zero creates no singularity at the chart boundary. If $z\in K_i$ and $w\notin\phi_i(U_i)$, then $\|z-w\|\geq\delta_i>0$, and the cross part of the seminorm is bounded by \[
2\int_{K_i}|(v_i^a\circ\phi_i^{-1})(z)|^p
\left(
\int_{\mathbb R^d\setminus\phi_i(U_i)}
\frac{dw}{\|z-w\|^{d+sp}}
\right)dz
\leq
C\delta_i^{-sp}
\|v_i^a\circ\phi_i^{-1}\|_{L^p(\phi_i(U_i))}^p.
\] The last inequality uses polar coordinates centered at $z$; the radial integral is $\displaystyle \int_{\delta_i}^{\infty}r^{-sp-1}\,dr=(sp)^{-1}\delta_i^{-sp}$. This proves the second inequality in \eqref{eq:caracterizacion-local-slobodeckij-intrinseco}.

For the reverse inequality, use $\displaystyle \mathbf{v}=\displaystyle\sum_{i=1}^{N}\psi_i\mathbf{v}$. The triangle inequality reduces the estimate to each section $\psi_i\mathbf{v}$. Pairs of points contained in $U_i$ are controlled by \eqref{eq:comparacion-local-slobodeckij-transporte-componentes}. If one point belongs to the support of $\psi_i$ and the other to $N\setminus U_i$, their distance is at least $\operatorname{dist}_{\mathbf{h}}(\operatorname{supp}\psi_i,N\setminus U_i)>0$; this part is controlled by the $L^p(N,\mathbf F)$ norm, as in Lemma~\ref{lem:independencia-radio-gagliardo-intrinseco}. This proves the first inequality when $0<s<1$.

Let $s=k+\sigma$, with $k\geq1$. In each chart, formula \eqref{eq:componentes-derivada-covariante-local}, applied to the localized section $\psi_i\mathbf{v}$, gives the triangular expression \begin{equation}
\label{eq:triangular-covariante-slobodeckij}
(\nabla_w^k(\psi_i\mathbf{v}))^a_{i_1\dots i_k}
=
\partial_{i_k}\cdots\partial_{i_1}v_i^a
+
\sum_{b=1}^r\sum_{|\beta|\leq k-1}
(B_{k,\beta})^a_{b\,i_1\dots i_k}\partial^\beta v_i^b.
\end{equation} The coefficients and all derivatives used are bounded on the fixed compact sets. Lemma~\ref{lem:W1p-controla-slobodeckij-global} controls each derivative of order at most $k-1$ in $W^{\sigma,p}(\mathbb R^d)$ by the $W^{k,p}(N,\mathbf F)$ norm, and its second assertion permits multiplication by the coefficients $B_{k,\beta}$ after inserting a cutoff equal to one on the support under consideration. Consequently, \eqref{eq:triangular-covariante-slobodeckij} and the formula obtained by solving for the partial derivatives of order $k$ give, in both directions, \[
\sum_{a=1}^r\|\widetilde v_i^a\|_{W^{k+\sigma,p}(\mathbb R^d)}
\leq C\|\psi_i\mathbf{v}\|_{W^{k+\sigma,p}(N,\mathbf{F})}
\] and \[
\|\psi_i\mathbf{v}\|_{W^{k+\sigma,p}(N,\mathbf{F})}
\leq C\sum_{a=1}^r
\|\widetilde v_i^a\|_{W^{k+\sigma,p}(\mathbb R^d)}.
\] Identity \eqref{eq:componentes-derivada-covariante-local} was established directly for weak covariant derivatives; thus \eqref{eq:triangular-covariante-slobodeckij} holds almost everywhere for the given section. Its inverse form is obtained by algebraically solving for the leading partial derivatives, without approximation in the fractional norm. Comparison of the fractional part of the components of $\nabla_w^k(\psi_i\mathbf{v})$ with its intrinsic seminorm follows by applying Lemma~\ref{lem:comparacion-local-slobodeckij-transporte-componentes} to the bundle $\mathbf{F}_k$ for pairs contained in $U_i$. The section $\nabla_w^k(\psi_i\mathbf{v})$ is supported in $\operatorname{supp}\psi_i\Subset U_i$; if only one of the two points is in $U_i$, the positive distance from the complement controls this part of the seminorm by the $L^p(N,\mathbf F_k)$ norm, exactly as in the case $0<s<1$. This gives the full intrinsic seminorm. Summing over the finite family of charts yields \eqref{eq:caracterizacion-local-slobodeckij-intrinseco} for every $k$. \end{proof}

\begin{corollary}[Functional structure and density] \label{cor:estructura-slobodeckij-intrinseco-secciones} Let $(N,\mathbf{h})$ be a compact connected Riemannian manifold without boundary, and let $\mathbf{F}\to N$ be a smooth vector bundle with a bundle metric and a compatible connection. Let $s>0$ and $1\leq p<\infty$. Then $W^{s,p}(N,\mathbf{F})$ is a Banach space. It is separable for $1\leq p<\infty$, reflexive for $1<p<\infty$, and $\Gamma(\mathbf{F})$ is dense in it. \end{corollary}

\begin{proof}
If $\dim N=0$, all the spaces under consideration are finite-dimensional and the assertions are immediate. Suppose that $d:=\dim N\geq1$.
 Proposition~\ref{prop:caracterizacion-local-slobodeckij-intrinseco} identifies the norm with a finite norm of Euclidean localizations. If $(\mathbf{v}_j)$ is Cauchy, each localized component converges in the Banach space $W^{s,p}(\mathbb R^d)$ by Theorem~\ref{teo:estructura-W-fraccionario}. The change-of-frame relations and the identity $\displaystyle \sum_{i=1}^{N}\psi_i=1$, valid for each $\mathbf{v}_j$, pass to the limit in $L^p(\mathbb R^d)$ and reconstruct a section $\mathbf{v}\in L^p(N,\mathbf{F})$. The first inequality in \eqref{eq:caracterizacion-local-slobodeckij-intrinseco} shows that $\mathbf{v}_j\to \mathbf{v}$ in $W^{s,p}(N,\mathbf{F})$.

For density, approximate each $\widetilde v_i^a$ in $W^{s,p}(\mathbb R^d)$ by functions in $C_c^\infty(\mathbb R^d)$ using the same theorem, and multiply the approximations by a fixed cutoff supported in $\phi_i(U_i)$ and equal to one on $K_i$. Reconstruction in the frame $\mathbf{e}_{i,a}$ and summation over $i,a$ give smooth sections; the localized estimate proves convergence. Separability follows from this same finite representation. For $1<p<\infty$, the integer-order case follows from
Corollary~\ref{Wmp uniformemente convexo}. Now let $s=k+\sigma$, with
$k\in\mathbb N_0$ and $0<\sigma<1$. Write
$\mathbf F_j=T^{(0,j)}(TN)\otimes\mathbf F$ for $0\leq j\leq k$, and set
\[
 Q_{\rho_N}:=\{(x,y)\in N\times N\mid 0<d_{\mathbf h}(x,y)<\rho_N\},
 \qquad \pi_1\colon Q_{\rho_N}\to N,\quad \pi_1(x,y)=x.
\]
Equip $Q_{\rho_N}$ with the restriction of the product measure
$\lambda_{\mathbf h}\otimes\lambda_{\mathbf h}$ and
$\pi_1^*\mathbf F_k$ with the induced metric. For a section $\mathbf w$ of
$\mathbf F_k$, define on this domain
\[
 (Q_\sigma\mathbf w)(x,y):=
 \frac{\mathbf w(x)-P^{\mathbf F_k}_{\gamma_{xy},1\to0}\mathbf w(y)}
 {d_{\mathbf h}(x,y)^{d/p+\sigma}}.
\]
The definition of the intrinsic norm shows that the linear map
\[
 \mathbf v\longmapsto
 \bigl((\nabla_w^j\mathbf v)_{j=0}^{k},Q_\sigma(\nabla_w^k\mathbf v)\bigr)
\]
is an isometric embedding into
\[
 \left(\prod_{j=0}^{k} L^p(N,\mathbf F_j)\right)
 \times L^p(Q_{\rho_N},\pi_1^*\mathbf F_k),
\]
equipped with the norm
\[
 \bigl\|((z_j)_{j=0}^k,z)\bigr\|
 :=\left(\sum_{j=0}^k\|z_j\|_{L^p(N,\mathbf F_j)}^p\right)^{1/p}
   +\|z\|_{L^p(Q_{\rho_N},\pi_1^*\mathbf F_k)}.
\]
The image is closed:
if the images of a sequence converge in this product, the isometry
makes the original sequence Cauchy in $W^{s,p}(N,\mathbf F)$;
the completeness already proved gives a limit in that space, and
continuity of the embedding identifies its image with the product limit.
Each factor is reflexive by uniform convexity of the spaces of integrable
sections when $1<p<\infty$; the Clarkson argument also applies to the
product measure and the bundle $\pi_1^*\mathbf F_k$.
Thus the finite product and its closed subspace are reflexive.
\end{proof}

\begin{remark}[Dependence on the connection] \label{obs:conexion-norma-slobodeckij-intrinseca} The seminorm in \eqref{eq:seminorma-gagliardo-intrinseca-secciones} records the connection and is therefore geometrically natural. Nevertheless, two smooth compatible connections on the same compact bundle yield the same space with equivalent norms. Indeed, Proposition~\ref{prop:caracterizacion-local-slobodeckij-intrinseco} compares both norms with the same Euclidean localizations; the constants may differ because they depend on the local connection matrices. \end{remark}

We now apply this construction to $N=\partial M$ and $\mathbf{F}=\mathbf{E}\restriction_{\partial M}$, with the induced metric and the connection from Lemma~\ref{lem:conexion-inducida-haz-restringido-frontera}. The resulting fractional regularity is therefore intrinsic to the geometric data of the boundary.

\begin{theorem}[Exact intrinsic trace for vector bundles] \label{teo:traza-slobodeckij-intrinseca-haces} \index{trace operator!in intrinsic Slobodeckij spaces} Let $(M,\mathbf{g})$ be a compact Riemannian manifold with nonempty smooth boundary, and let $\mathbf{E}\longrightarrow M$ be a smooth vector bundle of finite rank, equipped with a bundle metric $\mathbf{h}_{\mathbf{E}}$ and a compatible connection $\nabla^{\mathbf{E}}$. If $m\in\mathbb N$ and $1<p<\infty$, restriction of smooth sections extends uniquely to a continuous surjective linear operator \begin{equation}
\label{eq:traza-slobodeckij-intrinseca-haces}
\boldsymbol{\operatorname{Tr}}:
W^{m,p}(M,\mathbf{E})
\longrightarrow
W^{m-\frac1p,p}
(\partial M,\mathbf{E}\restriction_{\partial M}).
\end{equation} There is a continuous linear operator \begin{equation}
\label{eq:corretraccion-slobodeckij-intrinseca-haces}
\boldsymbol{\operatorname{Ex}}:
W^{m-\frac1p,p}
(\partial M,\mathbf{E}\restriction_{\partial M})
\longrightarrow W^{m,p}(M,\mathbf{E})
\end{equation} such that $\boldsymbol{\operatorname{Tr}}\boldsymbol{\operatorname{Ex}}=I$. In particular, there are constants $C_1,C_2>0$ such that \[
\|\boldsymbol{\operatorname{Tr}}\mathbf{u}\|_{W^{m-\frac1p,p}(\partial M,\mathbf{E}\restriction_{\partial M})}
\leq C_1\|\mathbf{u}\|_{W^{m,p}(M,\mathbf{E})}
\] and \[
\|\boldsymbol{\operatorname{Ex}}\mathbf{v}\|_{W^{m,p}(M,\mathbf{E})}
\leq C_2\|\mathbf{v}\|_{W^{m-\frac1p,p}(\partial M,\mathbf{E}\restriction_{\partial M})}.
\] \end{theorem}

\begin{proof} If $\dim M=1$, the boundary is a finite set and, under the preceding convention, the target space is the finite sum of the fibers of $\mathbf{E}$ over its points. On each interval component, the fundamental theorem of calculus and Hölder's inequality bound endpoint evaluation by the $W^{1,p}(M,\mathbf E)$ norm, and therefore by the $W^{m,p}(M,\mathbf E)$ norm. Choosing smooth sections supported in disjoint neighborhoods of the endpoints with prescribed values there gives a linear coretraction; its continuity is automatic because the target space is finite-dimensional. For the rest of the proof, assume that $\dim M\geq2$.

Set $s:=m-\frac{1}{p}=(m-1)+(1-\frac{1}{p})$; since $1<p<\infty$, its fractional part belongs to $(0,1)$. Choose the regular half-balls, frames, and finite partition $(\psi_\alpha)_{\alpha=1}^N$ from the proof of Theorem~\ref{teo:traza-entera-haces}, refined so that the tangential charts of $\partial M$ satisfy the hypotheses of Proposition~\ref{prop:caracterizacion-local-slobodeckij-intrinseco}.

First let $\mathbf{u}\in\Gamma(\mathbf{E})$. In a chart write \[
\psi_\alpha \mathbf{u}=\sum_{a=1}^rf_\alpha^a \mathbf{e}_{\alpha,a},
\qquad
v_\alpha^a:=f_\alpha^a\circ\phi_\alpha^{-1}
\quad\text{in }B_\alpha^+.
\] Each $v_\alpha^a$ vanishes on a neighborhood of the spherical part. Its extension by zero $\widetilde v_\alpha^a$ to $\mathbb R^n_+$ belongs to $W^{m,p}(\mathbb R^n_+)$ by Lemma~\ref{lem:extension-cero-media-bola-semiespacio}, and \[
\|\widetilde v_\alpha^a\|_{W^{m,p}(\mathbb R^n_+)}
=\|v_\alpha^a\|_{W^{m,p}(B_\alpha^+)}.
\] Apply the operator $\mathcal E_m$ from Corollary~\ref{teo:extension-sobolev-semiespacio}. The identification $W^{m,p}(\mathbb R^n)=H^{m,p}(\mathbb R^n)$ in Theorem~\ref{teo:identificaciones-H-F-W-entero}, followed by Theorem~\ref{teo:traza-fraccionaria-euclidea-hiperplano}, gives \[
\begin{aligned}
\|\operatorname{tr}(\mathcal E_m\widetilde v_\alpha^a)\|_
{B^{m-\frac1p}_{p,p}(\mathbb R^{n-1})}
&\leq C\|\mathcal E_m\widetilde v_\alpha^a\|_{H^{m,p}(\mathbb R^n)}\\
&\leq C'\|v_\alpha^a\|_{W^{m,p}(B_\alpha^+)}.
\end{aligned}
\] Since $m-\frac{1}{p}$ is not an integer, Corollary~\ref{teo:identificaciones-H-W-B-F} identifies $B^{m-\frac{1}{p}}_{p,p}(\mathbb R^{n-1})$ with $W^{m-\frac{1}{p},p}(\mathbb R^{n-1})$. On the flat ball, the preceding trace agrees with $(v_\alpha^a)\restriction_{B_\alpha'}$. Two extensions of $\widetilde v_\alpha^a$ agreeing on the half-space have the same trace: approximate in $W^{m,p}(\mathbb R^n_+)$ by functions smooth up to the boundary and use continuity of both trace operators; equivalently, apply the locality recorded in Lemma~\ref{lem:localidad-traza-semiespacio}. Since $\widetilde v_\alpha^a$ vanishes on a relative neighborhood of the complement of $B_\alpha'$ in the boundary of the half-space, the same lemma shows that the trace vanishes there. Thus $\operatorname{tr}(\mathcal E_m\widetilde v_\alpha^a)$ is exactly the extension by zero of the trace on $B_\alpha'$.

The localized characterization in Proposition~\ref{prop:caracterizacion-local-slobodeckij-intrinseco} now gives, with $s=m-\frac{1}{p}$, \[
\begin{aligned}
\|\mathbf{u}\restriction_{\partial M}\|_{W^{s,p}(\partial M,\mathbf{E}\restriction_{\partial M})}
&\leq C\sum_{\substack{1\leq\alpha\leq N\\1\leq a\leq r}}
\|\operatorname{tr}(\mathcal E_m\widetilde v_\alpha^a)\|_
{W^{s,p}(\mathbb R^{n-1})}\\
&\leq C\sum_{\substack{1\leq\alpha\leq N\\1\leq a\leq r}}
\|v_\alpha^a\|_{W^{m,p}(B_\alpha^+)}\\
&\leq C\sum_{\alpha=1}^{N}
\|\psi_\alpha \mathbf{u}\|_{W^{m,p}(U_\alpha,\mathbf{E})}\\
&\leq C\|\mathbf{u}\|_{W^{m,p}(M,\mathbf{E})}.
\end{aligned}
\] The third line uses local equivalence between partial and covariant derivatives; the fourth uses continuity of multiplication by each $\psi_\alpha$ and finiteness of the atlas. The constants can be chosen simultaneously because only finitely many charts, frames, and cutoffs are involved.

Theorem~\ref{teo:meyers-serrin-haz-frontera} and Corollary~\ref{cor:estructura-slobodeckij-intrinseco-secciones} allow restriction to extend uniquely and continuously to all of $W^{m,p}(M,\mathbf{E})$. The continuous inclusion \[
W^{m-\frac1p,p}(\partial M,\mathbf{E}\restriction_{\partial M})
\hookrightarrow W^{m-1,p}(\partial M,\mathbf{E}\restriction_{\partial M})
\] follows directly from the definition. By uniqueness of extension from smooth sections, the resulting operator agrees with that in Theorem~\ref{teo:traza-entera-haces}.

We construct the coretraction. Retain the preceding cover and choose a partition of unity $(\vartheta_\alpha)$ on $\partial M$ with $\operatorname{supp}\vartheta_\alpha\Subset U_\alpha\cap\partial M$. For $\mathbf{v}\in W^{m-\frac{1}{p},p}(\partial M,\mathbf{E}\restriction_{\partial M})$ write \[
\vartheta_\alpha \mathbf{v}
=\sum_{a=1}^rg_\alpha^a \mathbf{e}_{\alpha,a}
\] and extend $g_\alpha^a\circ(\phi_\alpha\restriction_{\partial M})^{-1}$ by zero to $\mathbb R^{n-1}$; denote this extension by $\widetilde g_\alpha^a$. Proposition \ref{prop:caracterizacion-local-slobodeckij-intrinseco} controls the sum of the $W^{m-\frac{1}{p},p}(\mathbb R^{n-1})$ norms of these extensions by the norm of $\mathbf{v}$.

Apply to each component the coretraction $\operatorname{ex}$ from Theorem~\ref{teo:traza-fraccionaria-euclidea-hiperplano}. The result belongs to $H^{m,p}(\mathbb R^n)=W^{m,p}(\mathbb R^n)$. Choose a cutoff $\chi_\alpha\in C_c^\infty(\mathbb R^n)$ whose support, intersected with the upper half-space, lies in the coordinate region, and whose restriction to the hyperplane equals one on a neighborhood of the supports of all localized components. Multiply $\operatorname{ex}\widetilde g_\alpha^a$ by $\chi_\alpha$, transport it to the chart, reconstruct it with the frame $\mathbf{e}_{\alpha,a}$, and extend it by zero outside $U_\alpha$. Denote the resulting global section by $\mathbf{V}_\alpha$. The support of the factor $\chi_\alpha\operatorname{ex}\widetilde g_\alpha^a$ is separated from the artificial part of the chart boundary; thus its transport to $U_\alpha$ and extension by zero have weak derivatives through order $m$ without boundary terms. This follows from the same compactly supported test argument as in Lemma~\ref{lem:extension-cero-media-bola-semiespacio}. The Euclidean multiplier bounds and local equivalence of norms give \[
\|\mathbf{V}_\alpha\|_{W^{m,p}(M,\mathbf{E})}
\leq C_\alpha\sum_{a=1}^r
\|\widetilde g_\alpha^a\|_{W^{m-\frac1p,p}(\mathbb R^{n-1})}.
\] Define \[
\boldsymbol{\operatorname{Ex}}\mathbf{v}
:=\sum_{\alpha=1}^{N} \mathbf{V}_\alpha.
\] The sum is finite, so the localized estimate proves continuity. Since $\chi_\alpha=1$ on the support of the boundary data and $\operatorname{tr}\operatorname{ex}=I$, we have \[
\boldsymbol{\operatorname{Tr}}\mathbf{V}_\alpha
=\vartheta_\alpha \mathbf{v}.
\] Summing over $\alpha$ gives $\displaystyle \boldsymbol{\operatorname{Tr}}\boldsymbol{\operatorname{Ex}}\mathbf{v}
=\displaystyle\sum_{\alpha=1}^{N}\vartheta_\alpha \mathbf{v}=\mathbf{v}$. This proves both surjectivity and the stated estimates. \end{proof} \begin{theorem}[Higher-order traces]\label{teo:trazas-orden-superior-haces}\index{higher-order traces} Let \((M,\mathbf{g})\) be a compact Riemannian manifold with nonempty boundary, and let \(\mathbf{E}\to M\) be a smooth vector bundle of finite rank, equipped with a bundle metric \(\mathbf{h}_{\mathbf{E}}\) and a compatible connection \(\nabla^{\mathbf{E}}\). Let \(m\geq1\) and \(1\leq p<\infty\). For each \(0\leq j\leq m-1\), the smooth restriction operator \[
\mathbf{u}\longmapsto \gamma_0(\nabla^j \mathbf{u})=(\nabla^j \mathbf{u})\restriction_{\partial M},
\qquad \mathbf{u}\in\Gamma(\mathbf{E}),
\] extends uniquely to a continuous linear operator \[
\boldsymbol{\operatorname{Tr}}_j\colon W^{m,p}(M,\mathbf{E})\longrightarrow
W^{m-j-1,p}\left(\partial M,\left(T^{(0,j)}(TM)\otimes \mathbf{E}\right)\restriction_{\partial M}\right).
\] Consequently, the operator \[
\mathcal T_m\colon W^{m,p}(M,\mathbf{E})\longrightarrow
\displaystyle\bigoplus_{j=0}^{m-1}
W^{m-j-1,p}\left(\partial M,\left(T^{(0,j)}(TM)\otimes \mathbf{E}\right)\restriction_{\partial M}\right),
\] defined by \(\mathcal T_m \mathbf{u}:=\bigl(\boldsymbol{\operatorname{Tr}}_0\mathbf{u},\boldsymbol{\operatorname{Tr}}_1\mathbf{u},\dots,\boldsymbol{\operatorname{Tr}}_{m-1}\mathbf{u}\bigr)\) is linear and continuous. If $1<p<\infty$, we also have the sharper estimate \[
\boldsymbol{\operatorname{Tr}}_j\colon
W^{m,p}(M,\mathbf{E})\longrightarrow
W^{m-j-\frac1p,p}
\left(\partial M,
\left(T^{(0,j)}(TM)\otimes \mathbf{E}\right)\restriction_{\partial M}
\right)
\] for every $j\in\{0,\dots,m-1\}$. \end{theorem}

\begin{proof} Fix \(0\leq j\leq m-1\) and write \(\mathbf{F}_j:=T^{(0,j)}(TM)\otimes \mathbf{E}\). This bundle is equipped with the product metric and induced connection. Moreover, \(\mathbf{F}_j\restriction_{\partial M}=\left(T^{(0,j)}(TM)\otimes \mathbf{E}\right)\restriction_{\partial M}\).

Let $\mathbf{u}\in W^{m,p}(M,\mathbf{E})$. By Theorem~\ref{teo:meyers-serrin-haz-frontera}, there is $(\mathbf{u}_\nu)\subseteq\Gamma(\mathbf{E})$ such that $\mathbf{u}_\nu\to \mathbf{u}$ in $W^{m,p}(M,\mathbf{E})$. For smooth sections, the recursive definition of covariant differentiation gives \[
\nabla^s(\nabla^j\mathbf{u}_\nu)=\nabla^{s+j}\mathbf{u}_\nu,
\qquad 0\leq s\leq m-j,
\] after the natural identification of cotangent factors. Consequently, $(\nabla^j\mathbf{u}_\nu)$ is Cauchy in $W^{m-j,p}(M,\mathbf{F}_j)$. Its limit in $L^p(M,\mathbf F_j)$ is $\nabla_w^j\mathbf{u}$; closedness of weak covariant derivatives then shows that $\nabla_w^j\mathbf{u}\in W^{m-j,p}(M,\mathbf{F}_j)$ and its successive derivatives are those expected. Moreover, \[
\|\nabla_w^j \mathbf{u}\|_{W^{m-j,p}(M,\mathbf{F}_j)}
\leq
\|\mathbf{u}\|_{W^{m,p}(M,\mathbf{E})}.
\]

Since \(m-j\geq1\), we may apply Theorem~\ref{teo:traza-entera-haces} to the bundle \(\mathbf{F}_j\to M\). This gives a continuous linear operator \[
\boldsymbol{\operatorname{Tr}}^{\mathbf{F}_j}:W^{m-j,p}(M,\mathbf{F}_j)
\longrightarrow
W^{m-j-1,p}(\partial M,\mathbf{F}_j\restriction_{\partial M}).
\] For \(\mathbf{u}\in W^{m,p}(M,\mathbf{E})\), define \(\boldsymbol{\operatorname{Tr}}_j \mathbf{u}:=\boldsymbol{\operatorname{Tr}}^{\mathbf{F}_j}(\nabla_w^j \mathbf{u})\). Then \[
\|\boldsymbol{\operatorname{Tr}}_j \mathbf{u}\|_{W^{m-j-1,p}(\partial M,\mathbf{F}_j\restriction_{\partial M})}
\leq
C\|\nabla_w^j \mathbf{u}\|_{W^{m-j,p}(M,\mathbf{F}_j)}
\leq
C\|\mathbf{u}\|_{W^{m,p}(M,\mathbf{E})}.
\] This proves that \(\boldsymbol{\operatorname{Tr}}_j\) is linear and continuous.

If $1<p<\infty$, instead apply Theorem~\ref{teo:traza-slobodeckij-intrinseca-haces} to the bundle $\mathbf{F}_j$ at order $m-j$. Since $\nabla_w^j\mathbf{u}\in W^{m-j,p}(M,\mathbf{F}_j)$, we obtain \[
\|\boldsymbol{\operatorname{Tr}}_j\mathbf{u}\|_
{W^{m-j-\frac1p,p}(\partial M,\mathbf{F}_j\restriction_{\partial M})}
\leq C\|\nabla_w^j\mathbf{u}\|_{W^{m-j,p}(M,\mathbf{F}_j)}
\leq C\|\mathbf{u}\|_{W^{m,p}(M,\mathbf{E})}.
\] After composition with the continuous inclusion into $W^{m-j-1,p}(\partial M,\mathbf{F}_j\restriction_{\partial M})$, the resulting operator and the one defined by the integer-order theorem agree on smooth sections and, by density, throughout $W^{m,p}(M,\mathbf{E})$.

We now show that this operator extends smooth restriction. If \(\mathbf{u}\in\Gamma(\mathbf{E})\), then \(\nabla_w^j \mathbf{u}=\nabla^j \mathbf{u}\), and by the definition of the trace operator in Theorem~\ref{teo:traza-entera-haces}, we have \(\boldsymbol{\operatorname{Tr}}^{\mathbf{F}_j}(\nabla^j \mathbf{u})=(\nabla^j \mathbf{u})\restriction_{\partial M}=\gamma_0(\nabla^j \mathbf{u})\). Thus \(\boldsymbol{\operatorname{Tr}}_j\) extends the smooth operator \(\mathbf{u}\mapsto\gamma_0(\nabla^j \mathbf{u})\).

Uniqueness follows from density of smooth sections in \(W^{m,p}(M,\mathbf{E})\). Finally, equip the direct sum \(\displaystyle\bigoplus_{j=0}^{m-1}W^{m-j-1,p}(\partial M,\mathbf{F}_j\restriction_{\partial M})\) with any of the equivalent product norms. Since each \(\boldsymbol{\operatorname{Tr}}_j\) is linear and continuous, there is a constant \(C_j>0\) such that \(\|\boldsymbol{\operatorname{Tr}}_j \mathbf{u}\|_{W^{m-j-1,p}(\partial M,\mathbf{F}_j\restriction_{\partial M})}\leq C_j\|\mathbf{u}\|_{W^{m,p}(M,\mathbf{E})}\). Summing over \(j\in\{0,\dots,m-1\}\) gives \[
\|\mathcal T_m \mathbf{u}\|
\leq
\left(\displaystyle\sum_{j=0}^{m-1}C_j\right)
\|\mathbf{u}\|_{W^{m,p}(M,\mathbf{E})}.
\] Therefore, \(\mathcal T_m\) is linear and continuous. \end{proof}

The integration-by-parts formula motivating this construction explains why \(\mathcal T_m\) is the appropriate boundary operator. Green's formula for \(\nabla^s\) contains terms of the form \[
\displaystyle\sum_{j=0}^{s-1}(-1)^{s-1-j}
\int_{\partial M}
\Big\langle
\nabla^j\mathbf{F},\,
\iota_{\boldsymbol{\nu}}\bigl((\operatorname{tr}_{\mathbf{g}}\circ\nabla)^{s-1-j}\mathbf{G}\bigr)
\Big\rangle_{\mathbf{g},\mathbf{h}_{\mathbf{E}}}
\,d\lambda_{\widetilde{\mathbf{g}}}.
\] Thus the quantities arising naturally on the boundary are restrictions of \(\nabla^j\mathbf{F}\) to \(\partial M\), with \(0\leq j\leq s-1\). This is why, to characterize \(W^{m,p}_0(M,\mathbf{E})\) when \(m\geq2\), vanishing of the trace of \(\mathbf{u}\) is insufficient: the traces of its covariant derivatives through order \(m-1\) must vanish as well.

\begin{theorem}[Characterization of \(W^{m,p}_0(M,\mathbf{E})\) by traces]\label{teo:caracterizacion-W0-trazas}\index{Sobolev space@Sobolev space!characterization of the zero-trace subspace}\index{trace operator!characterization of the kernel} Let \((M,\mathbf{g})\) be a compact Riemannian manifold with nonempty boundary, and let \(\mathbf{E}\to M\) be a smooth vector bundle of finite rank, equipped with a bundle metric \(\mathbf{h}_{\mathbf{E}}\) and a compatible connection \(\nabla^{\mathbf{E}}\). Let \(m\geq1\) and \(1\leq p<\infty\). Then \(W^{m,p}_0(M,\mathbf{E})=\ker(\mathcal T_m)\). That is, \[
W^{m,p}_0(M,\mathbf{E})
=
\left\{
\mathbf{u}\in W^{m,p}(M,\mathbf{E})
\middle|
\boldsymbol{\operatorname{Tr}}_j \mathbf{u}=0
\text{ for every }0\leq j\leq m-1
\right\}.
\] \end{theorem}

\begin{proof} We first prove that \(W^{m,p}_0(M,\mathbf{E})\subseteq\ker(\mathcal T_m)\). Let \(\mathbf{u}\in W^{m,p}_0(M,\mathbf{E})\). By definition, there is a sequence \((\mathbf{u}_k)_{k\in\mathbb N}\subseteq\Gamma_{c,\operatorname{Int}(M)}(\mathbf{E})\) such that \(\mathbf{u}_k\to \mathbf{u}\) in \(W^{m,p}(M,\mathbf{E})\). Since each \(\mathbf{u}_k\) has compact support contained in \(\operatorname{Int}(M)\), there is a neighborhood of \(\partial M\) on which \(\mathbf{u}_k\) vanishes. Thus, for every \(0\leq j\leq m-1\), \(\nabla^j \mathbf{u}_k\) also vanishes on a neighborhood of \(\partial M\). Consequently, \(\boldsymbol{\operatorname{Tr}}_j \mathbf{u}_k=0\). Using continuity of \(\boldsymbol{\operatorname{Tr}}_j\) gives \(\boldsymbol{\operatorname{Tr}}_j \mathbf{u}=\displaystyle\lim_{k\to\infty}\boldsymbol{\operatorname{Tr}}_j \mathbf{u}_k=0\). Since this holds for every \(0\leq j\leq m-1\), it follows that \(\mathbf{u}\in\ker(\mathcal T_m)\).

We now prove the reverse inclusion. Let \(\mathbf{u}\in\ker(\mathcal T_m)\). Take a finite open cover \((U_\alpha)_{\alpha=1}^{N}\) of \(M\) consisting of regular coordinate balls contained in \(\operatorname{Int}(M)\) and regular coordinate half-balls near \(\partial M\). Let \((\psi_\alpha)_{\alpha=1}^{N}\) be a smooth partition of unity subordinate to this cover, with \(\operatorname{supp}\psi_\alpha\Subset U_\alpha\). Since \(\mathbf{u}=\displaystyle\sum_{\alpha=1}^{N}\psi_\alpha \mathbf{u}\), it suffices to prove that \(\psi_\alpha \mathbf{u}\in W^{m,p}_0(M,\mathbf{E})\) for each \(\alpha\in\{1,\dots,N\}\).

Fix \(\alpha\in\{1,\dots,N\}\). First we prove that \(\boldsymbol{\operatorname{Tr}}_j(\psi_\alpha \mathbf{u})=0\) for every \(0\leq j\leq m-1\). Since \(\mathbf{u}\in\ker(\mathcal T_m)\), we have \(\boldsymbol{\operatorname{Tr}}_t \mathbf{u}=0\) for every \(0\leq t\leq m-1\). By Theorem~\ref{teo:meyers-serrin-haz-frontera}, there is a sequence \((\mathbf{v}_k)_{k\in\mathbb N}\subseteq\Gamma(\mathbf{E})\) such that \(\mathbf{v}_k\to \mathbf{u}\) in \(W^{m,p}(M,\mathbf{E})\). By continuity of the higher-order traces, given by Theorem~\ref{teo:trazas-orden-superior-haces}, we obtain \(\boldsymbol{\operatorname{Tr}}_t \mathbf{v}_k\to\boldsymbol{\operatorname{Tr}}_t \mathbf{u}=0\) in \(W^{m-t-1,p}\left(\partial M,\left(T^{(0,t)}(TM)\otimes \mathbf{E}\right)\restriction_{\partial M}\right)\) for every \(0\leq t\leq m-1\).

Now fix \(j\in\{0,\dots,m-1\}\). For smooth sections, the iterated Leibniz rule for the product connection gives, up to the natural permutations of covariant indices, \[
\nabla^j(\psi_\alpha \mathbf{v}_k)
=
\displaystyle\sum_{t=0}^{j}
\binom{j}{t}
(\nabla^{j-t}\psi_\alpha)\otimes\nabla^t \mathbf{v}_k.
\] More precisely, this equality means that the right-hand side represents a finite sum of terms obtained from \((\nabla^{j-t}\psi_\alpha)\otimes\nabla^t \mathbf{v}_k\) by permutations of the \(j\) covariant indices. Restricting to \(\partial M\) gives \[
\boldsymbol{\operatorname{Tr}}_j(\psi_\alpha \mathbf{v}_k)
=
\displaystyle\sum_{t=0}^{j}
\binom{j}{t}
\bigl((\nabla^{j-t}\psi_\alpha)\restriction_{\partial M}\bigr)\otimes \boldsymbol{\operatorname{Tr}}_t \mathbf{v}_k,
\] again up to those index permutations. This notation suffices for the estimate we need, because index permutations are pointwise isometries, and tensor multiplication by the fixed smooth sections \((\nabla^{j-t}\psi_\alpha)\restriction_{\partial M}\) is continuous on the corresponding Sobolev spaces over the compact manifold \(\partial M\). Since \(\boldsymbol{\operatorname{Tr}}_t \mathbf{v}_k\to0\) for every \(0\leq t\leq j\), we conclude that \(\boldsymbol{\operatorname{Tr}}_j(\psi_\alpha \mathbf{v}_k)\to0\) in \(W^{m-j-1,p}\left(\partial M,\left(T^{(0,j)}(TM)\otimes \mathbf{E}\right)\restriction_{\partial M}\right)\).

On the other hand, the map \(\mathbf{v}\longmapsto\boldsymbol{\operatorname{Tr}}_j(\psi_\alpha \mathbf{v})\) is continuous from \(W^{m,p}(M,\mathbf{E})\) to that same space: this follows from continuity of multiplication by \(\psi_\alpha\) on \(W^{m,p}(M,\mathbf{E})\), given by Lemma~\ref{lem:multiplicacion-funcion-suave-sobolev}, and continuity of \(\boldsymbol{\operatorname{Tr}}_j\). Since \(\mathbf{v}_k\to \mathbf{u}\) in \(W^{m,p}(M,\mathbf{E})\), we have \(\boldsymbol{\operatorname{Tr}}_j(\psi_\alpha \mathbf{v}_k)\to\boldsymbol{\operatorname{Tr}}_j(\psi_\alpha \mathbf{u})\). By uniqueness of the limit, \(\boldsymbol{\operatorname{Tr}}_j(\psi_\alpha \mathbf{u})=0\). Since \(j\) was arbitrary, we obtain \(\boldsymbol{\operatorname{Tr}}_j(\psi_\alpha \mathbf{u})=0\) for every \(0\leq j\leq m-1\).

Now consider the two types of charts. If \(U_\alpha\) is a regular coordinate ball contained in \(\operatorname{Int}(M)\), then \(\psi_\alpha \mathbf{u}\) has compact support contained in \(\operatorname{Int}(M)\). Take a function \(\chi_\alpha\in C_c^\infty(\operatorname{Int}(M))\) such that \(\chi_\alpha\equiv1\) on a neighborhood of \(\operatorname{supp}(\psi_\alpha)\). By Theorem~\ref{teo:meyers-serrin-haz-frontera}, there is a sequence \((\mathbf{z}_k)_{k\in\mathbb N}\subseteq\Gamma(\mathbf{E})\) such that \(\mathbf{z}_k\to\psi_\alpha \mathbf{u}\) in \(W^{m,p}(M,\mathbf{E})\). By Lemma~\ref{lem:multiplicacion-funcion-suave-sobolev}, we have \(\chi_\alpha \mathbf{z}_k\to\chi_\alpha\psi_\alpha \mathbf{u}=\psi_\alpha \mathbf{u}\) in \(W^{m,p}(M,\mathbf{E})\). Moreover, each \(\chi_\alpha \mathbf{z}_k\) is a smooth section with compact support contained in \(\operatorname{Int}(M)\). Therefore, \(\psi_\alpha \mathbf{u}\in W^{m,p}_0(M,\mathbf{E})\).

Now suppose that \(U_\alpha\) is a regular coordinate half-ball. Write \(B_\alpha^+:=\phi_\alpha(U_\alpha\cap\operatorname{Int}(M))\) and \(B_\alpha':=\phi_\alpha(U_\alpha\cap\partial M)\). Take a smooth local frame \((\mathbf{e}_{1,\alpha},\dots,\mathbf{e}_{r,\alpha})\) of \(\mathbf{E}\) on \(U_\alpha\), and write locally \(\psi_\alpha \mathbf{u}=\displaystyle\sum_{a=1}^{r}w_\alpha^a\mathbf{e}_{a,\alpha}\). Set \(q_\alpha^a:=w_\alpha^a\circ\phi_\alpha^{-1}\) on \(B_\alpha^+\), and denote by \(\widetilde q_\alpha^a\) its extension by zero to \(\mathbb R^n_+\). Since \(\operatorname{supp}\psi_\alpha\Subset U_\alpha\), each \(q_\alpha^a\) vanishes outside a compact subset of \(B_\alpha^+\cup B_\alpha'\) lying at positive distance from the spherical part of \(B_\alpha^+\). Thus we may use the preceding convention for these functions.

We want to prove that each \(q_\alpha^a\) belongs to the closure of \(C_c^\infty(B_\alpha^+)\) in \(W^{m,p}(B_\alpha^+)\). By Corollary~\ref{cor:caracterizacion-local-media-bola}, it suffices to prove that \(\gamma_0D^\beta q_\alpha^a=0\) on \(B_\alpha'\), in the sense of the preceding convention, for every \(|\beta|\leq m-1\) and every \(a\). Explicitly, we must prove that \((\gamma_0D^\beta\widetilde q_\alpha^a)\restriction_{B_\alpha'}=0\).

We prove this by induction on \(\ell\). For \(\ell=0\), the equality \(\boldsymbol{\operatorname{Tr}}_0(\psi_\alpha \mathbf{u})=0\) means that the trace of \(\psi_\alpha \mathbf{u}\) vanishes as a section of \(\mathbf{E}\restriction_{\partial M}\). In the boundary chart and chosen frame, the trace operator was constructed precisely by taking components, extending them by zero to the half-space, and applying the half-space trace. Therefore, \((\gamma_0\widetilde q_\alpha^a)\restriction_{B_\alpha'}=0\) for every \(a\).

Suppose that \((\gamma_0D^\beta\widetilde q_\alpha^a)\restriction_{B_\alpha'}=0\) has been proved for every \(|\beta|\leq\ell-1\) and every \(a\), where \(1\leq\ell\leq m-1\). Since \(\boldsymbol{\operatorname{Tr}}_\ell(\psi_\alpha \mathbf{u})=0\), the local components of the trace of \(\nabla^\ell(\psi_\alpha \mathbf{u})\) vanish on \(B_\alpha'\). By Lemma~\ref{lem:local-expression-higher-order}, applied in the chart \((U_\alpha,\phi_\alpha)\) and frame \((\mathbf{e}_{1,\alpha},\dots,\mathbf{e}_{r,\alpha})\), we have locally \[
(\nabla^\ell(\psi_\alpha \mathbf{u}))^a_{i_1\dots i_\ell}
=
\partial_{i_\ell}\cdots\partial_{i_1}w_\alpha^a
+
\displaystyle\sum_{b=1}^{r}
\displaystyle\sum_{|\beta|\leq\ell-1}
(B_\beta)^a_{b\,i_1\dots i_\ell}\,\partial^\beta w_\alpha^b.
\] To pass this identity to traces, first fix the supports of the approximation. Choose $\chi_\alpha\in C_c^\infty(U_\alpha)$ equal to one near $\operatorname{supp}\psi_\alpha$ and take $\mathbf z_k\in\Gamma(\mathbf E)$ with $\mathbf z_k\to\psi_\alpha\mathbf u$ in $W^{m,p}(M,\mathbf E)$, by Theorem~\ref{teo:meyers-serrin-haz-frontera}. Then $\mathbf v_k=\chi_\alpha\mathbf z_k$ converges to the same limit in that space by Lemma~\ref{lem:multiplicacion-funcion-suave-sobolev}, and all its supports lie in the fixed compact set $\operatorname{supp}\chi_\alpha$. Denote its coordinate components by $q_{\alpha,k}^a$. Local equivalence of norms and Lemma~\ref{lem:extension-cero-media-bola-semiespacio} give \[
 \widetilde q_{\alpha,k}^a\longrightarrow\widetilde q_\alpha^a
 \quad\text{in }W^{m,p}(\mathbb R^n_+).
\] Here extension by zero is from the half-ball to the half-space, across its spherical part; the fixed compact set is separated from that part, and the flat boundary is not crossed.

For $|\beta|\leq\ell$, $D^\beta\widetilde q_{\alpha,k}^a\to
D^\beta\widetilde q_\alpha^a$ in $W^{m-|\beta|,p}(\mathbb R^n_+)$, and by continuity of the trace, \[
 \gamma_0D^\beta\widetilde q_{\alpha,k}^a
 \longrightarrow\gamma_0D^\beta\widetilde q_\alpha^a
 \quad\text{in }W^{m-|\beta|-1,p}(\mathbb R^{n-1}).
\] All these spaces embed continuously into $W^{m-\ell-1,p}(\mathbb R^{n-1})$. Multiply the coefficients $B_\beta$, expressed in the chart, by a cutoff equal to one near the fixed compact set; they then extend to smooth compactly supported coefficients on $\mathbb R^n$. Multiplication by them, before and after taking the trace, is continuous. Moreover, the components of $\nabla^\ell\mathbf v_k$ converge in $W^{m-\ell,p}(\mathbb R^n_+)$, since the covariant derivative of order $\ell$ is continuous from $W^{m,p}(M,\mathbf E)$ to $W^{m-\ell,p}(M,T^{(0,\ell)}(TM)\otimes\mathbf E)$, as proved in Theorem~\ref{teo:trazas-orden-superior-haces}. Apply the local identity to $\mathbf v_k$ and take limits in $W^{m-\ell-1,p}(\mathbb R^{n-1})$, then restrict to $B_\alpha'$. This gives exactly the required trace identity for $\psi_\alpha\mathbf u$, including when $p=1$.

Now fix a multi-index \(\gamma\) with \(|\gamma|=\ell\). Choose an ordered block \((i_1,\dots,i_\ell)\) in which the index \(q\) appears exactly \(\gamma_q\) times. Since weak partial derivatives commute, the leading term of the corresponding component is \[
\partial_{i_\ell}\cdots\partial_{i_1}q_\alpha^a
=D^\gamma q_\alpha^a.
\] Taking the trace of that component on \(B_\alpha'\) gives \[
0
=
\gamma_0\!\left(D^\gamma\widetilde q_\alpha^a\right)
+
\sum_{b=1}^{r}\sum_{|\beta|\leq\ell-1}
\gamma_0\!\left(
\widetilde B_{\beta}\,D^\beta\widetilde q_\alpha^b
\right)
\quad\text{in }B_\alpha',
\] where \(\widetilde B_\beta\) denotes the smooth coefficient expressed in the chart. The left-hand side is zero because it is a component of \(\boldsymbol{\operatorname{Tr}}_\ell(\psi_\alpha\mathbf u)\). For each term in the sum, compatibility of the trace with multiplication by a smooth function gives \[
\gamma_0\!\left(
\widetilde B_{\beta}\,D^\beta\widetilde q_\alpha^b
\right)
=
\widetilde B_{\beta}\restriction_{B_\alpha'}\,
\gamma_0\!\left(D^\beta\widetilde q_\alpha^b\right)
=0,
\] since \(|\beta|\leq\ell-1\) and the induction hypothesis applies. This compatibility uses the same sequence with fixed support: if $v_k\to v$ in $W^{m-|\beta|,p}(\mathbb R^n_+)$, both $\gamma_0(\widetilde B_\beta v_k)$ and $(\widetilde B_\beta|_{\mathbb R^{n-1}})\gamma_0v_k$ converge in $W^{m-|\beta|-1,p}(\mathbb R^{n-1})$ by the multiplication and trace bounds. They agree for each smooth $k$; uniqueness of the limit proves the identity used. Consequently, \[
\gamma_0\left(D^\gamma\widetilde q_\alpha^a\right)\restriction_{B_\alpha'}
=
0
\] for every \(|\gamma|=\ell\) and every \(a\). This completes the induction.

Thus \(\gamma_0D^\beta q_\alpha^a=0\) on \(B_\alpha'\), in the sense of the preceding convention, for every \(|\beta|\leq m-1\) and every \(a\). By Corollary~\ref{cor:caracterizacion-local-media-bola}, each \(q_\alpha^a\) belongs to the closure of \(C_c^\infty(B_\alpha^+)\) in \(W^{m,p}(B_\alpha^+)\). Hence, for each \(a\) there is a sequence \((\varphi_k^a)_{k\in\mathbb N}\subseteq C_c^\infty(B_\alpha^+)\) such that \(\varphi_k^a\to q_\alpha^a\) in \(W^{m,p}(B_\alpha^+)\).

On \(U_\alpha\cap\operatorname{Int}(M)\) define \(\mathbf{s}_k:=\displaystyle\sum_{a=1}^{r}(\varphi_k^a\circ\phi_\alpha)\mathbf{e}_{a,\alpha}\). Since each \(\varphi_k^a\) has compact support contained in \(B_\alpha^+\), each \(\mathbf{s}_k\) extends by zero to a smooth section with compact support contained in \(\operatorname{Int}(M)\). By local equivalence of norms in the regular coordinate half-ball, we have \(\mathbf{s}_k\to\psi_\alpha \mathbf{u}\) in \(W^{m,p}(M,\mathbf{E})\). Therefore, \(\psi_\alpha \mathbf{u}\in W^{m,p}_0(M,\mathbf{E})\).

We have proved that \(\psi_\alpha \mathbf{u}\in W^{m,p}_0(M,\mathbf{E})\) for each \(\alpha\in\{1,\dots,N\}\). Since \(\mathbf{u}=\displaystyle\sum_{\alpha=1}^{N}\psi_\alpha \mathbf{u}\), we conclude that \(\mathbf{u}\in W^{m,p}_0(M,\mathbf{E})\). This proves that \(\ker(\mathcal T_m)\subseteq W^{m,p}_0(M,\mathbf{E})\) and completes the proof. \end{proof}

\subsection{Embeddings for sections on compact manifolds with boundary}

The local estimates will be applied to components supported in a
compact subset of the chart domain, even when it meets the flat face.
Lemma~\ref{lem:esquinas-encajes-localizados-ortante} supplies the
scalar estimate and compactness with this fixed support; extension
of sections allows us to combine them globally.

\begin{theorem}[Sobolev embedding for bundles on manifolds with boundary]
\label{teo:encaje-sobolev-haces-frontera}
\index{teorema de encaje de sobolev@Sobolev embedding theorem!for bundles on manifolds with boundary}
Let $(M,\mathbf g)$ be a compact Riemannian manifold of dimension
$n$ with nonempty smooth boundary, and let $\mathbf E\to M$ be a
finite rank real or complex bundle with smooth fiber metric and
compatible connection. Let $j\geq0$, $m\geq1$ be integers and
$1\leq p,q<\infty$ satisfy
\[
 j+m-\frac{n}{p}\geq j-\frac{n}{q}.
\]
Then
\[
 W^{j+m,p}(M,\mathbf E)\hookrightarrow W^{j,q}(M,\mathbf E)
\]
continuously, without boundary conditions.
\end{theorem}

\begin{proof}
Choose finitely many regular charts $(U_\alpha,\phi_\alpha)$,
frames of rank $h$ on their larger domains, and cutoffs
$\psi_\alpha$ compactly supported in $U_\alpha$, using
Proposition~\ref{prop:esquinas-cartas-regulares-cubierta}.
Let $1\leq\alpha\leq N$, $r_\alpha\in\{0,1\}$, and
\[
 \Omega_\alpha:=\phi_\alpha(U_\alpha)\subseteq G_{n,r_\alpha},
 \quad\Omega_{\alpha,+}:=\phi_\alpha(U_\alpha\cap\operatorname{Int}(M)),
 \quad K_\alpha:=\phi_\alpha(\operatorname{supp}\psi_\alpha).
\]
Once these data are fixed, none of the compact sets depend on
the section being estimated.

For $\mathbf u\in W^{j+m,p}(M,\mathbf E)$, apply the Leibniz rule
of Lemma~\ref{lem:esquinas-cortes-extension-cero}
to each local component. The derivatives of
$\psi_\alpha\circ\phi_\alpha^{-1}$ are bounded on the chart.
Both directions of Lemma~\ref{lema:meyers-serrin-haz-frontera}
control the product in the source covariant norm.
Its support is compact in $U_\alpha$, so
Lemma~\ref{lem:frontera-extension-cero-secciones-soporte}
extends it to the global space with the same norm. We obtain
\[
 \|\psi_\alpha\mathbf u\|_{W^{j+m,p}(M,\mathbf E)}
                 \leq D_\alpha\|\mathbf u\|_{W^{j+m,p}(M,\mathbf E)}.
\]
If $\psi_\alpha\mathbf u=\sum_a u_\alpha^a\mathbf e_{\alpha,a}$,
put $v_\alpha^a=u_\alpha^a\circ\phi_\alpha^{-1}$.
The local characterization at order $j+m$ and exponent $p$
shows that these functions belong to $W^{j+m,p}(\Omega_{\alpha,+})$
and
\[
 c_\alpha\sum_{a=1}^h
       \|v_\alpha^a\|_{W^{j+m,p}(\Omega_{\alpha,+})}
 \leq\|\psi_\alpha\mathbf u\|_{W^{j+m,p}(U_\alpha,\mathbf E|_{U_\alpha})}.
\]
Each component vanishes outside $K_\alpha$.
Apply Lemma~\ref{lem:esquinas-encajes-localizados-ortante} with
$\ell=j+m$, both for $r_\alpha=0$ and for $r_\alpha=1$:
\[
 v_\alpha^a\in W^{j,q}(\Omega_{\alpha,+}),\qquad
 \|v_\alpha^a\|_{W^{j,q}(\Omega_{\alpha,+})}
       \leq A_\alpha\|v_\alpha^a\|_{W^{j+m,p}(\Omega_{\alpha,+})}.
\]
The reverse local characterization, now at order $j$ and exponent
$q$, gives
$(\psi_\alpha\mathbf u)|_{U_\alpha}\in W^{j,q}(U_\alpha,\mathbf E|_{U_\alpha})$
and
\[
 \|\psi_\alpha\mathbf u\|_{W^{j,q}(U_\alpha,\mathbf E|_{U_\alpha})}
       \leq C_\alpha\sum_{a=1}^h
                 \|v_\alpha^a\|_{W^{j,q}(\Omega_{\alpha,+})}.
\]
Its essential support lies in
$\operatorname{supp}\psi_\alpha\Subset U_\alpha$.
Lemma~\ref{lem:frontera-extension-cero-secciones-soporte}
extends it to the global target space, preserving the norm.
This extension is exactly $\psi_\alpha\mathbf u$.
Combining the preceding bounds yields
\[
 \|\psi_\alpha\mathbf u\|_{W^{j,q}(M,\mathbf E)}
 \leq\frac{C_\alpha A_\alpha D_\alpha}{c_\alpha}
                      \|\mathbf u\|_{W^{j+m,p}(M,\mathbf E)}.
\]
Since $\mathbf u=\sum_{\alpha=1}^N\psi_\alpha\mathbf u$ and every
summand already belongs to the target space, so does $\mathbf u$.
The triangle inequality gives
\[
 \|\mathbf u\|_{W^{j,q}(M,\mathbf E)}
 \leq\left(\sum_{\alpha=1}^N
       \frac{C_\alpha A_\alpha D_\alpha}{c_\alpha}\right)
                         \|\mathbf u\|_{W^{j+m,p}(M,\mathbf E)}.
\]
The finite sum is a constant independent of $\mathbf u$.
\end{proof}

\begin{theorem}[Rellich--Kondrashov for bundles on manifolds with boundary]
\label{teo:rellich-kondrashov-haces-frontera}
\index{teorema de rellich kondrashov@Rellich--Kondrashov theorem!for bundles on manifolds with boundary}
Under the same geometric hypotheses, let $j\geq0$, $m\geq1$
be integers and $1\leq p,q<\infty$. If
\[
 j+m-\frac{n}{p}>j-\frac{n}{q},
\]
the embedding
$W^{j+m,p}(M,\mathbf E)\hookrightarrow W^{j,q}(M,\mathbf E)$
is compact.
\end{theorem}

\begin{proof}
Let $(\mathbf u_k)$ be bounded in the source space, with
$\|\mathbf u_k\|_{W^{j+m,p}(M,\mathbf E)}\leq C_0$.
Use the fixed charts, cutoffs, and frames of the preceding proof,
and put
\[
 v_{\alpha,k}^a:=
       (\psi_\alpha\mathbf u_k)_\alpha^a\circ\phi_\alpha^{-1}.
\]
The multiplier and local characterization estimates give
\[
 \sum_{a=1}^h\|v_{\alpha,k}^a\|_{W^{j+m,p}(\Omega_{\alpha,+})}
                    \leq c_\alpha^{-1}D_\alpha C_0.
\]
All supports lie in the same $K_\alpha$, even if it meets the
flat face. The compactness assertion of
Lemma~\ref{lem:esquinas-encajes-localizados-ortante} gives a
subsequence converging in $W^{j,q}(\Omega_{\alpha,+})$ for
each pair $(\alpha,a)$. Enumerate the $Nh$ pairs and successively
extract a subsequence of the previous one.
After these finitely many extractions, one subsequence is Cauchy
for all pairs.

The embedding theorem already proved places every $\mathbf u_k$
in $W^{j,q}(M,\mathbf E)$. For $\varepsilon>0$, choose a common
index $k_0$ such that, for $k,l\geq k_0$ and every $\alpha$,
\[
 \sum_{a=1}^h\|v_{\alpha,k}^a-v_{\alpha,l}^a
                    \|_{W^{j,q}(\Omega_{\alpha,+})}
                 <\frac{\varepsilon}{NC_\alpha}.
\]
The localized difference is supported in
$\operatorname{supp}\psi_\alpha$, and the extension lemma
identifies its local and global norms.
Norm equivalence in the target space implies
\[
 \begin{aligned}
 \|\mathbf u_k-\mathbf u_l\|_{W^{j,q}(M,\mathbf E)}
 &\leq\sum_{\alpha=1}^N
             \|\psi_\alpha(\mathbf u_k-\mathbf u_l)\|_{W^{j,q}(M,\mathbf E)}\\
 &=\sum_{\alpha=1}^N
             \|\psi_\alpha(\mathbf u_k-\mathbf u_l)\|_{W^{j,q}(U_\alpha,\mathbf E|_{U_\alpha})}\\
 &\leq\sum_{\alpha=1}^NC_\alpha\sum_{a=1}^h
       \|v_{\alpha,k}^a-v_{\alpha,l}^a\|_{W^{j,q}(\Omega_{\alpha,+})}
 <\varepsilon.
 \end{aligned}
\]
The subsequence is Cauchy in a Banach space and converges.
This proves compactness.
\end{proof}

\begin{corollary}[Embeddings and the critical exponent] \label{cor:encajes-criticos-haces-frontera} Under the geometric hypotheses of the preceding theorems, let $j\geq0$ and $m\geq1$ be integers, and let $1\leq p<n/m$. Then \[
W^{j+m,p}(M,\mathbf E)\hookrightarrow W^{j,q}(M,\mathbf E)
\] continuously for $1\leq q\leq p_m^*:=\displaystyle\frac{np}{n-mp}$, and compactly for $1\leq q<p_m^*$. \end{corollary} \begin{proof} Since $mp<n$, the condition $j+m-\displaystyle\frac np\geq j-\displaystyle\frac nq$ is equivalent to $q\leq p_m^*$, and the strict inequality is equivalent to $q<p_m^*$. Apply Theorems~\ref{teo:encaje-sobolev-haces-frontera} and \ref{teo:rellich-kondrashov-haces-frontera}, respectively. \end{proof}

The scalar case follows by taking the trivial rank-one bundle with its canonical metric and connection.

\begin{corollary}[Sobolev embedding on compact manifolds with boundary] \label{teo:encaje-sobolev-variedades-frontera} \index{Sobolev embedding theorem@Sobolev embedding theorem!on manifolds with boundary} Let $(M,\mathbf{g})$ be a compact Riemannian manifold of dimension $n$ with nonempty smooth boundary. Let $j\geq0$ and $m\geq1$ be integers and $1\leq p,q<\infty$ such that $j+m-\displaystyle\frac np\geq j-\displaystyle\frac nq$. Then the inclusion \[
W^{j+m,p}(M)\hookrightarrow W^{j,q}(M)
\] is continuous. \end{corollary}

\begin{proof} Apply Theorem~\ref{teo:encaje-sobolev-haces-frontera} to the trivial bundle $M\times\mathbb R\to M$, equipped with the usual metric and the connection induced by the exterior differential. Under the canonical identification of its sections with functions, the Sobolev norms agree with the scalar norms. \end{proof}

\begin{corollary}[Rellich--Kondrashov on compact manifolds with boundary] \label{teo:rellich-kondrashov-variedades-frontera} \index{Rellich Kondrashov theorem@Rellich--Kondrashov theorem!on manifolds with boundary} Let $(M,\mathbf{g})$ be a compact Riemannian manifold of dimension $n$ with nonempty smooth boundary. Let $j\geq0$ and $m\geq1$ be integers and $1\leq p,q<\infty$ such that $j+m-\frac{n}{p}>j-\frac{n}{q}$. Then the inclusion \[
W^{j+m,p}(M)\hookrightarrow W^{j,q}(M)
\] is compact. \end{corollary}

\begin{proof} This is the case of Theorem~\ref{teo:rellich-kondrashov-haces-frontera} corresponding to the trivial rank-one bundle. \end{proof}

\begin{corollary}[Compact embedding with a decrease of one order] \label{cor:rellich-descenso-un-orden-haces} Let $(M,\mathbf{g})$ be a compact Riemannian manifold with possibly empty smooth boundary, and let $\mathbf{E}\to M$ be a smooth real or complex vector bundle of finite rank, equipped with a bundle metric $\mathbf{h}_{\mathbf{E}}$ (Hermitian in the complex case) and a compatible connection $\nabla^{\mathbf{E}}$. For every integer $m\geq1$ and every $1\leq p<\infty$, the inclusions \[
W^{m,p}(M,\mathbf{E})\hookrightarrow W^{m-1,p}(M,\mathbf{E}),
\qquad
W_0^{m,p}(M,\mathbf{E})\hookrightarrow W_0^{m-1,p}(M,\mathbf{E})
\] are compact. Here $W_0^{0,p}(M,\mathbf{E})=L^p(M,\mathbf{E})$. \end{corollary} \begin{proof} In Theorems~\ref{teo: rellich--kondrashov para haces vectoriales} and \ref{teo:rellich-kondrashov-haces-frontera}, take $j=m-1$, a gain of one derivative, and the same exponent $p$ in the source and target. The inequality between indices reduces to $1>0$, giving the first inclusion for every $1\leq p<\infty$, including $p=1$.

For the second, if $\mathbf{u}\in W_0^{m,p}(M,\mathbf{E})$, choose $\mathbf{u}_\nu\in\Gamma_{c,\operatorname{Int}(M)}(\mathbf{E})$ with $\mathbf{u}_\nu\to\mathbf{u}$ in $W^{m,p}(M,\mathbf{E})$. The inequality $\|\mathbf{v}\|_{W^{m-1,p}(M,\mathbf E)}\leq\|\mathbf{v}\|_{W^{m,p}(M,\mathbf E)}$ implies that the same convergence holds in $W^{m-1,p}(M,\mathbf{E})$; by definition, $\mathbf{u}\in W_0^{m-1,p}(M,\mathbf{E})$. By the first inclusion, a bounded sequence in $W_0^{m,p}(M,\mathbf E)$ has a subsequence converging in $W^{m-1,p}(M,\mathbf E)$. All its terms belong to $W_0^{m-1,p}(M,\mathbf E)$, and this subspace is closed by definition, so the limit also belongs to it. The subspace norm is the induced norm, proving the asserted compactness. \end{proof}

\begin{remark}[Total orders and dependence of the constants] \label{obs:rellich-general-haces-compactos} In terms of total orders, the two preceding theorems state that \[
W^{m,p}(M,\mathbf{E})\hookrightarrow W^{k,q}(M,\mathbf{E})
\quad\text{is compact if}\quad
m>k\geq0,\quad 1\leq p,q<\infty,\quad
m-\frac np>k-\frac nq.
\] This formulation is obtained by substituting $j=k$ and replacing the gain of derivatives by $m-k$. The same assertion holds with subscript zero in both source and target: continuity of the inclusion carries smooth approximations supported in the interior into $W^{k,q}(M,\mathbf E)$, and closedness of $W_0^{k,q}(M,\mathbf E)$ retains the limit. If $(m-k)p<n$, the strict condition is equivalent to $q<np/(n-(m-k)p)$; at the critical exponent, only continuity has been asserted in the corresponding Sobolev theorems.

The continuity constants depend on $m,k,p,q$, $M$, $\mathbf{g}$, and the smooth data $\mathbf{E}$, $\mathbf{h}_{\mathbf{E}}$, $\nabla^{\mathbf{E}}$. In the proofs they arise from a fixed finite cover, derivatives of its partition functions, and bounds on the coefficients of the metrics, their inverses, and the connections through the orders used. They do not depend on the section or the sequence index. For a varying family of geometric data, uniformity of these constants requires control of those bounds and the charts; compactness of each manifold separately does not provide a uniform constant for the entire family. No assertion in this paragraph includes infinite exponents; embeddings into continuous and Hölder sections are treated separately. \end{remark}

\begin{corollary}[First-order embeddings for the Dirichlet space] \label{cor:encajes-W0-haces-frontera} Let $(M,\mathbf{g})$ be a compact Riemannian manifold of dimension $n$ with nonempty smooth boundary, and let $\mathbf{E}\to M$ be a smooth real vector bundle of finite rank, equipped with a bundle metric $\mathbf{h}_{\mathbf{E}}$ and a connection $\nabla^{\mathbf{E}}$ compatible with $\mathbf{h}_{\mathbf{E}}$. Let $1\leq q<n$ and define $\displaystyle q^*:=\displaystyle\frac{nq}{n-q}$. Then \[
W_0^{1,q}(M,\mathbf{E})\hookrightarrow L^\ell(M,\mathbf{E})
\] continuously for every $1\leq\ell\leq q^*$, and compactly for every $1\leq\ell<q^*$. \end{corollary}

\begin{proof} The continuous inclusion follows from Theorem~\ref{teo:encaje-sobolev-haces-frontera} with $j=0$ and $m=1$. The compact inclusion follows from Theorem~\ref{teo:rellich-kondrashov-haces-frontera}. In both cases it suffices to restrict the inclusion operators to the closed subspace $W_0^{1,q}(M,\mathbf{E})\subseteq W^{1,q}(M,\mathbf{E})$. \end{proof}

To formulate Morrey's embedding up to the boundary, we need to extend the notation of Definition~\ref{def:secciones-holder-haz}. We do not use the seminorm defined by minimizing geodesics, since that formulation assumes the manifold has no boundary. Instead, we use the Euclidean extension norms already fixed for $C^{k,\alpha}(\overline\Omega)$.

\begin{definition}[Hölder sections up to the boundary] \label{def:secciones-holder-haz-frontera} \index{Holder section@Hölder section!up to the boundary} Let $(M,\mathbf{g})$ be a compact Riemannian manifold with nonempty smooth boundary, and let $\mathbf{E}\to M$ be a smooth vector bundle of rank $r$, equipped with a bundle metric and a compatible connection. Choose a finite cover $(U_\lambda,\phi_\lambda)_{\lambda=1}^N$ consisting of regular interior coordinate balls and regular boundary coordinate half-balls, a smooth partition of unity $(\psi_\lambda)_{\lambda=1}^N$ subordinate to it, with \(
\operatorname{supp}\psi_\lambda\Subset U_\lambda
\), and local frames $(\mathbf{e}_{\lambda,1},\dots,\mathbf{e}_{\lambda,r})$ defined on neighborhoods of $\overline{U_\lambda}$. Set \[
\Omega_\lambda
:=
\begin{cases}
\phi_\lambda(U_\lambda),&U_\lambda\subseteq\operatorname{Int}(M),\\
\phi_\lambda(U_\lambda\cap\operatorname{Int}(M)),
&U_\lambda\cap\partial M\neq\varnothing.
\end{cases}
\] If $\mathbf{u}\in\Gamma^k(\mathbf{E})$ and \[
\psi_\lambda \mathbf{u}
=
\displaystyle\sum_{a=1}^r u_\lambda^a \mathbf{e}_{\lambda,a}
\qquad\text{in }U_\lambda,
\] define \begin{equation}
\label{eq:norma-holder-localizada-haces-frontera}
\|\mathbf{u}\|_{\Gamma^{k,\alpha}_{\mathcal A}(\mathbf{E})}
:=
\displaystyle\sum_{\lambda=1}^N\displaystyle\sum_{a=1}^r
\|u_\lambda^a\circ\phi_\lambda^{-1}\|_{C^{k,\alpha}(\overline{\Omega_\lambda})},
\end{equation} where $k\in\mathbb N_0$, $0<\alpha<1$, and $\mathcal A$ represents the chosen set of local data. Denote by $\Gamma^{k,\alpha}_{\mathcal A}(\mathbf{E})$ the space of sections of class $C^k$ up to the boundary for which \eqref{eq:norma-holder-localizada-haces-frontera} is finite. \end{definition}

\begin{proposition}[Independence of the Hölder norm up to the boundary] \label{prop:independencia-norma-holder-localizada-haces-frontera} Let $(M,\mathbf{g})$ be a compact Riemannian manifold with nonempty smooth boundary, and let $\mathbf{E}\to M$ be a smooth vector bundle of finite rank with a bundle metric and a compatible connection. Let $k\in\mathbb N_0$ and $0<\alpha<1$. The space and topology in Definition~\ref{def:secciones-holder-haz-frontera} are independent of the chosen cover, partition of unity, charts, and frames. Moreover, the resulting space is Banach. Henceforth we denote it by $\Gamma^{k,\alpha}(\mathbf{E})$ and omit the subscript $\mathcal A$ from its norm. \end{proposition}

\begin{proof} Let $\mathcal A$ and $\widetilde{\mathcal A}$ be two choices of local data. Use a tilde for objects in the second system. For an index $\mu$ of $\widetilde{\mathcal A}$, the identity $\displaystyle\sum_{\lambda=1}^N\psi_\lambda=1$ gives \begin{equation}
\label{eq:transferencia-holder-frontera-particion}
\widetilde\psi_\mu \mathbf{u}
=
\sum_{\substack{1\leq\lambda\leq N\\U_\lambda\cap\widetilde U_\mu\neq\varnothing}}
\widetilde\psi_\mu\psi_\lambda \mathbf{u}.
\end{equation} The sum is finite. On an overlap write \[
\mathbf{e}_{\lambda,a}
=
\sum_{b=1}^r
(G_{\mu\lambda})^b_{\ a}\,\widetilde{\mathbf{e}}_{\mu,b}
\] and denote the coordinate change on its domain by \[
\Phi_{\lambda\mu}
:=
\phi_\lambda\circ\widetilde\phi_\mu^{-1}
\]. If $f_\lambda^a:=u_\lambda^a\circ\phi_\lambda^{-1}$, the component $b$ of the corresponding summand in \eqref{eq:transferencia-holder-frontera-particion}, transported to the chart $\widetilde\phi_\mu$, is \begin{equation}
\label{eq:operador-transferencia-holder-frontera}
\left[
(\widetilde\psi_\mu G_{\mu\lambda})^b_{\ a}
\circ\widetilde\phi_\mu^{-1}
\right]
\bigl(f_\lambda^a\circ\Phi_{\lambda\mu}\bigr).
\end{equation}

We explain precisely why this expression defines a continuous operator between the extension norms. The support of the first factor lies in a compact subset of the domain of the coordinate change. In the chart $\widetilde\phi_\mu$, choose a smooth cutoff equal to one on a neighborhood of that compact set, with support still contained in the same domain. Changes of charts on a manifold with boundary extend smoothly to Euclidean open sets containing those compact sets; for two boundary charts the extension crosses the flat face, while in an overlap with an interior chart the compact set lies away from that face. Take an extension $F_\lambda^a\in C^{k,\alpha}(\mathbb R^n)$ of $f_\lambda^a$. After inserting the preceding cutoff, \eqref{eq:operador-transferencia-holder-frontera} is the restriction of a product of a smooth compactly supported function and $F_\lambda^a$ composed with a smooth diffeomorphism.

Denote by $T_{\mu\lambda}^{ba}f_\lambda^a$ the function on the right-hand side of \eqref{eq:operador-transferencia-holder-frontera}. The Faà di Bruno formula in Theorem~\ref{faa di bruno multivariable} controls derivatives of the composition through order $k$. For derivatives of order $k$, we also use \[
[fg]_{C^{0,\alpha}(K)}
\leq
\|f\|_{C^0(K)}[g]_{C^{0,\alpha}(K)}
+
\|g\|_{C^0(K)}[f]_{C^{0,\alpha}(K)}
\] and the fact that a smooth map is Lipschitz on a compact set. We obtain a constant $C_{\lambda\mu}>0$ such that \[
\|T_{\mu\lambda}^{ba}f_\lambda^a\|_{C^{k,\alpha}(\overline{\widetilde\Omega_\mu})}
\leq
C_{\lambda\mu}
\|f_\lambda^a\|_{C^{k,\alpha}(\overline{\Omega_\lambda})}.
\] The inequality is first obtained for an extension $F_\lambda^a$, and then the infimum over all extensions is taken in the norm definition.

Sum over the finitely many overlaps, components, and indices. The maximum of the preceding constants is finite, and \eqref{eq:transferencia-holder-frontera-particion} gives \[
\|\mathbf{u}\|_{\Gamma^{k,\alpha}_{\widetilde{\mathcal A}}(\mathbf{E})}
\leq
C\|\mathbf{u}\|_{\Gamma^{k,\alpha}_{\mathcal A}(\mathbf{E})}.
\] Interchanging the two systems gives the reverse inequality.

We prove completeness. If $(\mathbf{u}_j)_{j\in\mathbb N}$ is Cauchy for one of these norms, each family of localized components converges in the Banach space $C^{k,\alpha}(\overline{\Omega_\lambda})$, whose completeness was established in Proposition~\ref{prop:completitud-holder-euclidiano}. Reconstruct the limits in each frame, multiply them by cutoffs equal to one on the supports of $\psi_\lambda$, and extend them by zero outside $U_\lambda$. The finite sum of the resulting sections defines a section $\mathbf{u}$ of class $C^k$ up to the boundary. Since $\mathbf{u}_j=\displaystyle\sum_{\lambda=1}^N\psi_\lambda \mathbf{u}_j$, the same transfer estimate shows that $\mathbf{u}_j\to \mathbf{u}$ in $\Gamma^{k,\alpha}(\mathbf{E})$. \end{proof}

\begin{theorem}[Sobolev--Morrey embedding for bundles over manifolds with boundary] \label{teo:sobolev-morrey-haces-frontera} \index{Sobolev Morrey theorem@Sobolev--Morrey theorem!for bundles over manifolds with boundary} Let $(M,\mathbf{g})$ be a compact Riemannian manifold of dimension $n$ with nonempty smooth boundary, and let $\mathbf{E}\to M$ be a smooth vector bundle of finite rank, equipped with a bundle metric and a compatible connection. Let $m\in\mathbb N$, $1\leq p<\infty$, $k\in\mathbb N_0$, and $0<\alpha<1$. If \begin{equation}
\label{eq:condicion-sobolev-morrey-haces-frontera}
k+\alpha
\leq
m-\frac np,
\end{equation} then each element of $W^{m,p}(M,\mathbf{E})$ has a unique representative in $\Gamma^{k,\alpha}(\mathbf{E})$, and there is $C>0$ such that \begin{equation}
\label{eq:estimacion-sobolev-morrey-haces-frontera}
\|\mathbf{u}\|_{\Gamma^{k,\alpha}(\mathbf{E})}
\leq
C\|\mathbf{u}\|_{W^{m,p}(M,\mathbf{E})}.
\end{equation} The restriction of this representative to $\partial M$ agrees with $\boldsymbol{\operatorname{Tr}}\mathbf{u}$. \end{theorem}

\begin{proof} Use the local data from Definition~\ref{def:secciones-holder-haz-frontera}. For $\mathbf{u}\in W^{m,p}(M,\mathbf{E})$ write \[
\psi_\lambda \mathbf{u}
=
\sum_{a=1}^r u_\lambda^a \mathbf{e}_{\lambda,a}
\qquad\text{in }U_\lambda
\] and set $v_\lambda^a:=u_\lambda^a\circ\phi_\lambda^{-1}$. Continuity of multiplication by $\psi_\lambda$, local equivalence between ordinary and covariant derivatives, and finiteness of the family give a constant $C_0>0$, independent of $\mathbf{u}$, such that \begin{equation}
\label{eq:cota-componentes-morrey-haces-frontera}
\sum_{\lambda=1}^N\sum_{a=1}^r
\|v_\lambda^a\|_{W^{m,p}(\Omega_\lambda)}
\leq
C_0\|\mathbf{u}\|_{W^{m,p}(M,\mathbf{E})}.
\end{equation}

If $U_\lambda$ is an interior chart, the support of $v_\lambda^a$ is compact in $\Omega_\lambda$ and therefore $v_\lambda^a\in W^{m,p}_0(\Omega_\lambda)$. Theorem~\ref{teo:sobolev-morrey-euclidiano} gives \begin{equation}
\label{eq:morrey-local-interior-haz-frontera}
\|v_\lambda^a\|_{C^{k,\alpha}(\overline{\Omega_\lambda})}
\leq
C_\lambda
\|v_\lambda^a\|_{W^{m,p}(\Omega_\lambda)}.
\end{equation}

Suppose that $U_\lambda$ is a boundary chart and write $\Omega_\lambda=B_\lambda^+$. There is a compact set \[
K_\lambda
\subseteq
B_\lambda^+\cup B_\lambda'
\] containing the essential supports of all components $v_\lambda^a$ and lying at positive distance from the spherical part of $\partial B_\lambda^+$. Lemma~\ref{lem:extension-cero-media-bola-semiespacio} permits extension of $v_\lambda^a$ by zero across that spherical part, but not across the flat face. Denote this first extension by $w_\lambda^a$. Then \[
w_\lambda^a\in W^{m,p}(\mathbb R^n_+),
\qquad
\|w_\lambda^a\|_{W^{m,p}(\mathbb R^n_+)}
=
\|v_\lambda^a\|_{W^{m,p}(B_\lambda^+)}.
\] Now apply the extension operator \[
\mathcal E_m\colon W^{m,p}(\mathbb R^n_+)\longrightarrow W^{m,p}(\mathbb R^n)
\] from Corollary~\ref{teo:extension-sobolev-semiespacio}. Choose $R_\lambda>0$ and $\chi_\lambda\in C_c^\infty(B_{\mathrm{euc}}(0,2R_\lambda))$ such that $\chi_\lambda=1$ on a neighborhood of $K_\lambda$, and define \[
T_\lambda v_\lambda^a
:=
\chi_\lambda\mathcal E_m w_\lambda^a.
\] This function belongs to $W^{m,p}_0(B_{\mathrm{euc}}(0,2R_\lambda))$, its restriction to $B_\lambda^+$ agrees with $v_\lambda^a$, and \[
\|T_\lambda v_\lambda^a\|_{W^{m,p}(B_{\mathrm{euc}}(0,2R_\lambda))}
\leq
C_\lambda'
\|v_\lambda^a\|_{W^{m,p}(B_\lambda^+)}.
\] Theorem~\ref{teo:sobolev-morrey-euclidiano}, applied to $T_\lambda v_\lambda^a$, and the definition of the extension norm give \begin{equation}
\label{eq:morrey-local-frontera-haz}
\|v_\lambda^a\|_{C^{k,\alpha}(\overline{B_\lambda^+})}
\leq
C_\lambda''
\|v_\lambda^a\|_{W^{m,p}(B_\lambda^+)}.
\end{equation}

For each $\lambda$, the continuous representatives of $v_\lambda^1,\dots,v_\lambda^r$ reconstruct a localized section $\mathbf{s}_\lambda$ on $U_\lambda$. This section agrees almost everywhere with $\psi_\lambda \mathbf{u}$ and, since its components vanish outside the fixed compact set determined by $\operatorname{supp}(\psi_\lambda)$, it extends by zero to a global section of class $C^k$ up to the boundary. Define the global representative by \[
\widetilde{\mathbf{u}}:=\sum_{\lambda=1}^N \mathbf{s}_\lambda.
\] The sum is finite and agrees almost everywhere with $\displaystyle\sum_{\lambda=1}^N\psi_\lambda \mathbf{u}=\mathbf{u}$. Definition~\ref{def:secciones-holder-haz-frontera}, the local estimates, and \eqref{eq:cota-componentes-morrey-haces-frontera} give \[
\|\widetilde{\mathbf{u}}\|_{\Gamma^{k,\alpha}(\mathbf{E})}
\leq
C\sum_{\lambda=1}^N\sum_{a=1}^r
\|v_\lambda^a\|_{W^{m,p}(\Omega_\lambda)}
\leq
C'\|\mathbf{u}\|_{W^{m,p}(M,\mathbf{E})}.
\] If two representatives continuous up to the boundary define the same class in $L^p(M,\mathbf{E})$, the relative set where they differ is open. Every nonempty relative coordinate ball or half-ball has positive Riemannian volume, so that set is empty. This proves uniqueness.

It remains to identify restriction to the boundary. By Theorem~\ref{teo:meyers-serrin-haz-frontera}, choose $(\mathbf{u}_j)_{j\in\mathbb N}\subseteq\Gamma(\mathbf{E})$ such that $\mathbf{u}_j\to \mathbf{u}$ in $W^{m,p}(M,\mathbf{E})$. Applying \eqref{eq:estimacion-sobolev-morrey-haces-frontera} to differences shows that $\mathbf{u}_j$ converges to the preceding representative in $\Gamma^{k,\alpha}(\mathbf{E})$ and, in particular, uniformly on $\partial M$. On the other hand, continuity in Theorem~\ref{teo:traza-entera-haces} implies \[
\mathbf{u}_j\restriction_{\partial M}
=
\boldsymbol{\operatorname{Tr}}\mathbf{u}_j
\longrightarrow
\boldsymbol{\operatorname{Tr}}\mathbf{u}
\qquad\text{in }L^p(\partial M,\mathbf{E}\restriction_{\partial M}).
\] Uniform convergence gives the same limit in $L^p(\partial M,\mathbf E\restriction_{\partial M})$; uniqueness of the limit identifies the restriction of the representative with the trace. \end{proof}

\begin{theorem}[Rellich--Kondrashov into Hölder sections up to the boundary] \label{teo:rellich-kondrashov-holder-haces-frontera} \index{Rellich Kondrashov theorem@Rellich--Kondrashov theorem!into Hölder sections up to the boundary} Let $(M,\mathbf{g})$ be a compact Riemannian manifold of dimension $n$ with nonempty smooth boundary, and let $\mathbf{E}\to M$ be a smooth vector bundle of finite rank with a bundle metric and a compatible connection. Let $m\in\mathbb N$, $1\leq p<\infty$, $k\in\mathbb N_0$, and $0<\alpha<1$. If \[
k+\alpha
<
m-\frac np,
\] then the embedding \[
W^{m,p}(M,\mathbf{E})
\hookrightarrow
\Gamma^{k,\alpha}(\mathbf{E})
\] is compact. \end{theorem}

\begin{proof} Let $(\mathbf{u}_j)_{j\in\mathbb N}$ be a bounded sequence in $W^{m,p}(M,\mathbf{E})$, and use the localizations from the preceding proof. In an interior chart, each component sequence is bounded in $W^{m,p}_0(\Omega_\lambda)$; Theorem~\ref{teo:rellich-kondrashov-holder-euclidiano} gives a subsequence converging in $C^{k,\alpha}(\overline{\Omega_\lambda})$.

In a boundary chart, apply to each component the operator $T_\lambda$ constructed in the preceding proof. The functions $T_\lambda v_{\lambda,j}^a$ form a bounded sequence in $W^{m,p}_0(B_{\mathrm{euc}}(0,2R_\lambda))$. The same Euclidean theorem gives a subsequence converging in $C^{k,\alpha}(\overline{B_{\mathrm{euc}}(0,2R_\lambda)})$; upon restriction, the original components converge in $C^{k,\alpha}(\overline{B_\lambda^+})$.

There are only finitely many charts and components. A diagonal extraction produces a single subsequence for which all localized components are Cauchy in their Hölder spaces. By \eqref{eq:norma-holder-localizada-haces-frontera}, the subsequence is Cauchy in $\Gamma^{k,\alpha}(\mathbf{E})$ and converges by completeness established in Proposition~\ref{prop:independencia-norma-holder-localizada-haces-frontera}. \end{proof}

\begin{corollary}[Regularity exponents up to the boundary] \label{cor:exponentes-regulares-holder-haces-frontera} Let $(M,\mathbf{g})$ be a compact Riemannian manifold of dimension $n$ with nonempty smooth boundary, and let $\mathbf{E}\to M$ be a smooth vector bundle of finite rank with a bundle metric and a compatible connection. Let $m\in\mathbb N$ and $1\leq p<\infty$, and let $\sigma:=m-\frac{n}{p}>0$. \begin{enumerate}[label=(\alph*)] \item If $\sigma=\ell+\theta$, with $\ell\in\mathbb N_0$ and $0<\theta<1$, then \[
W^{m,p}(M,\mathbf{E})
\hookrightarrow
\Gamma^{\ell,\theta}(\mathbf{E})
\] continuously, and the embedding into $\Gamma^{\ell,\beta}(\mathbf{E})$ is compact for every $0<\beta<\theta$. \item If $\sigma=s\in\mathbb N$, then, for every $0<\alpha<1$, \[
W^{m,p}(M,\mathbf{E})
\hookrightarrow
\Gamma^{s-1,\alpha}(\mathbf{E})
\] compactly and, in particular, continuously. \end{enumerate} \end{corollary} \begin{proof} For \(\sigma=\ell+\theta\), the Sobolev--Morrey theorem up to the boundary, Theorem~\ref{teo:sobolev-morrey-haces-frontera}, applied with \(k=\ell\) and \(\alpha=\theta\), gives the continuous embedding. If \(0<\beta<\theta\), then \[
\ell+\beta<m-\frac np,
\] and Theorem~\ref{teo:rellich-kondrashov-holder-haces-frontera} gives the compact embedding into \(\Gamma^{\ell,\beta}(\mathbf E)\).

If \(\sigma=s\in\mathbb N\) and \(0<\alpha<1\), then \[
(s-1)+\alpha<s=m-\frac np.
\] Applying the Rellich--Kondrashov theorem up to the boundary once more with \(k=s-1\) gives the compact, and hence continuous, embedding into \(\Gamma^{s-1,\alpha}(\mathbf E)\). \end{proof}

\section{Sobolev spaces on manifolds with corners}
\label{sec:sobolev-haces-esquinas}

The product in Theorem~\ref{teo:producto-variedades-con-esquinas}
requires models with several normal coordinates. We use the
ordinary tangent bundle, of full dimension, constructed in
Proposition~\ref{prop:esquinas-tangente-conos}.

\begin{proposition}[Geometric data, measure, and adjoints at corners]
\label{prop:esquinas-geometria-medida-adjunto}
Every smooth manifold with corners admits smooth positive
Riemannian metrics on its ordinary tangent bundle. Its finite rank
vector bundles admit smooth fiber metrics and compatible
connections. For a fixed metric $\mathbf g$, the Riemann--Lebesgue
measure $\lambda_{\mathbf g}$ is defined in charts by the density
$\sqrt{\det(\mathbf g)}$, and
$\lambda_{\mathbf g}(\partial M)=0$.
The Levi--Civita connection and tensor connections are smooth up
to every face. Every differential operator with coefficients smooth
up to the faces has a unique interior formal adjoint, smooth up
to the faces, with respect to the fixed metrics.
\end{proposition}

\begin{proof}
Patch local Euclidean metrics using the partition in
Proposition~\ref{prop:esquinas-particiones-unidad}. A sum of
positive forms with nonnegative weights adding to one is positive
on every fiber. The same procedure constructs fiber metrics,
Hermitian for complex bundles. Gram--Schmidt gives smooth local
orthonormal frames. In each frame, take the flat connection and
patch them using the partition: the Leibniz rule and metric
compatibility are preserved because the weights add to one.

In a chart, points of positive depth lie in a finite union of
coordinate hyperplanes, of zero $n$-dimensional Lebesgue measure.
A countable cover shows that the boundary is null.
In the interior, change of variables patches the volume densities;
adding the null faces defines the same measure independently of
the charts. The coefficients
\[
 \Gamma_{ij}^k=\frac{1}{2}\sum_{\ell=1}^n g^{k\ell}
       \bigl(\partial_i g_{j\ell}+\partial_j g_{i\ell}
                                    -\partial_\ell g_{ij}\bigr)
\]
are smooth up to the faces: the metric matrix is invertible, and
coordinate derivatives are unique by
Lemma~\ref{lem:esquinas-jets-extension-local}.
In the interior, they give the Levi--Civita connection. The coordinate
transformation laws and the torsion and compatibility identities
extend to the faces by continuity and density of the interior.
The same holds for product and dual connections.

Finally, let $H_{\mathbf E},H_{\mathbf F}$ be the matrices of two
fiber metrics, $\omega=\sqrt{\det(\mathbf g)}$, and
$P=\sum_\beta A_\beta D^\beta$ the coordinate expression of an
operator. Integration against compactly supported interior tests gives
\[
 P_h^*\Psi=(\omega H_{\mathbf E})^{-1}
       \sum_\beta(-1)^{|\beta|}
           D^\beta(A_\beta^*\omega H_{\mathbf F}\Psi).
\]
All its coefficients have smooth local extensions, and the inverses
exist near every point, including points on the faces.
On a chart overlap, both expressions satisfy the adjunction identity
for interior tests. Uniqueness there follows by choosing tests with
arbitrary jets at a point and using uniqueness of distributions
represented by smooth functions. Hence the expressions agree in
the interior. Their smooth coefficients also agree on the faces
by continuity. They therefore patch to a global operator, and
the same argument proves uniqueness.
\end{proof}

\begin{definition}[Weak derivatives and norms at corners]
\label{def:sobolev-haces-esquinas}
Let $(M,\mathbf g)$ be a smooth Riemannian manifold with corners
of dimension $n$,
and let $\mathbf E\to M$ be a finite rank real or complex vector
bundle with fiber metric $\mathbf h_{\mathbf E}$ and connection
smooth up to the faces. The iterates $\nabla^s$ use the Levi--Civita
and product connections, placing the new covariant slot first.
For $\mathbf u\in L^1_{\operatorname{loc}}(M,\mathbf E)$,
a section $\mathbf v_s$ is $\nabla_w^s\mathbf u$ if
\[
 \int_M\langle\mathbf v_s,\boldsymbol\sigma\rangle_{\mathbf g,\mathbf h_{\mathbf E}}
                              \,d\lambda_{\mathbf g}
 =\int_M\langle\mathbf u,(\nabla^s)_h^*\boldsymbol\sigma\rangle_{\mathbf h_{\mathbf E}}
                              \,d\lambda_{\mathbf g}
\]
for every
$\boldsymbol\sigma\in\Gamma_c(\operatorname{Int}(M),
T^{(0,s)}(TM)\otimes\mathbf E)$. Set
$\nabla_w^0\mathbf u=\mathbf u$. For $m\in\mathbb N_0$ and
$1\leq p\leq\infty$, require these derivatives to exist in $L^p$
for $0\leq s\leq m$, and define
\[
 \|\mathbf u\|_{W^{m,p}(M,\mathbf E)}
 :=\begin{cases}
 \displaystyle\left(\sum_{s=0}^m\int_M
       |\nabla_w^s\mathbf u|_{\mathbf g,\mathbf h_{\mathbf E}}^p
                         \,d\lambda_{\mathbf g}\right)^{\frac{1}{p}},
                                      &p<\infty,\\
 \displaystyle\max_{0\leq s\leq m}
       \|\nabla_w^s\mathbf u\|_{L^\infty(M,T^{(0,s)}(TM)\otimes\mathbf E)},
                                      &p=\infty.
 \end{cases}
\]
Local regularity is required on precompact neighborhoods that
may meet every face. Equivalence classes and norms do not depend
on the values assigned on $\partial M$.
\end{definition}

\begin{lemma}[Local characterization in a regular corner chart]
\label{lem:esquinas-componentes-normas-locales}
Let $(U,\phi)$ be a regular chart modeled on $G_{n,r}$, and put
$\Omega:=\phi(U)$, $\Omega_+:=\Omega\cap O_{n,r}$.
Assume that the frame of $\mathbf E$ and the geometric data are
defined on a neighborhood of $\overline U$. If
$\mathbf u=\sum_{a=1}^h u^a\mathbf e_a$ and
$f^a=u^a\circ\phi^{-1}$, then, for
$m\in\mathbb N_0$ and $1\leq p\leq\infty$,
\[
 \mathbf u\in W^{m,p}(U,\mathbf E|_U)
 \quad\Longleftrightarrow\quad
 f^a\in W^{m,p}(\Omega_+)\quad(1\leq a\leq h),
\]
and
\[
 c\sum_{a=1}^h\|f^a\|_{W^{m,p}(\Omega_+)}
 \leq\|\mathbf u\|_{W^{m,p}(U,\mathbf E|_U)}
 \leq C\sum_{a=1}^h\|f^a\|_{W^{m,p}(\Omega_+)}.
\]
The constants depend on the chart and data, but not on $\mathbf u$.
\end{lemma}

\begin{proof}
Depth in this chart identifies
$\phi(U\cap\operatorname{Int}(M))=\Omega_+$.
The product connection components satisfy
\[
 V_{s+1,iI}^a
 =\partial_iV_{s,I}^a+\sum_b A_{i,b}^aV_{s,I}^b
   -\sum_{\nu=1}^s\sum_k\Gamma_{ii_\nu}^kV_{s,I(\nu:k)}^a.
\]
Induction, starting with $V_0^a=f^a$, gives
\[
 (P_sf)_I^a=D_If^a+
       \sum_{b=1}^h\sum_{|\beta|<s}B_{s,I,\beta,b}^aD^\beta f^b.
\]
Indeed, in passing from $s$ to $s+1$, only differentiation of
$D_If^a$ produces order $s+1$; all remaining terms have order at
most $s$. Thus the coefficients are constructed by sums, products,
and derivatives of smooth coefficients.

For an interior tensor test, Proposition~
\ref{prop:esquinas-geometria-medida-adjunto} gives
$(P_s)_h^*\Psi=(\omega H_0)^{-1}P_s^\dagger(\omega H_s\Psi)$.
The substitution $\zeta=\omega H_s\Psi$ is a bijection on the
compactly supported smooth tests on $\Omega_+$.
The weak identity is therefore equivalent to
\[
 \int_{\Omega_+}\zeta^*V_s\,dx
      =\int_{\Omega_+}(P_s^\dagger\zeta)^*f\,dx,
 \qquad V_s=P_sf.
\]
The second equality is understood in the distributional sense.
The tests stay away from all faces, regardless of the number $r$.
If $f$ has $L^p$ partial derivatives, the expression for $P_s$
represents every weak covariant derivative. Conversely, the identity
\[
 D_If^a=(V_s)_I^a-
       \sum_{b=1}^h\sum_{|\beta|<s}B_{s,I,\beta,b}^aD^\beta f^b
\]
recovers, by induction on $s$, every partial derivative as an
$L^p$ function. All multi-indices of length $s$ occur among the blocks $I$.

The chart closure is compact in the larger domain. On it, the
coefficients are bounded, the volume densities have positive bounds,
and the metric matrices have uniformly positive and bounded
eigenvalues. The direct expression for $P_s$ bounds its components
by partial derivatives; the solved expression and induction bound
partial derivatives by covariant components of orders at most $s$.
Taking $L^p$ norms and combining the finite families proves both
inequalities, as in the two explicit estimates of
Lemma~\ref{lema:meyers-serrin-haz-frontera}. For $p=\infty$, use
finite maxima. The nullity of the faces identifies interior norms
with the norms on $U$.
\end{proof}

\begin{lemma}[Locality, multiplication, and extension of sections]
\label{lem:esquinas-localidad-cortes-secciones}
For the spaces in Definition~\ref{def:sobolev-haces-esquinas},
restriction to an open set preserves weak derivatives.
These derivatives vanish almost everywhere wherever the section
vanishes on an open set. For a smooth function $\chi$ whose covariant
derivatives through order $m$ are bounded,
\[
 |\nabla_w^s(\chi\mathbf u)|
 \leq\sum_{t=0}^s\binom{s}{t}
            |\nabla^{s-t}\chi|\,|\nabla_w^t\mathbf u|,
 \qquad 0\leq s\leq m,
\]
and
\[
 \|\chi\mathbf u\|_{W^{m,p}(M,\mathbf E)}
 \leq C_{m,p}\max_{0\leq t\leq m}
              \|\nabla^t\chi\|_{L^\infty(M)}
              \|\mathbf u\|_{W^{m,p}(M,\mathbf E)}.
\]
In particular, this applies to every $\chi\in C_c^\infty(M)$.
If $\mathbf v\in W^{m,p}(U,\mathbf E|_U)$ vanishes outside a
compact set $K\subseteq U$, its zero extension satisfies
\[
 \nabla_w^s(Z_U\mathbf v)=Z_U(\nabla_w^s\mathbf v),\qquad
 \|Z_U\mathbf v\|_{W^{m,p}(M,\mathbf E)}
       =\|\mathbf v\|_{W^{m,p}(U,\mathbf E|_U)}.
\]
\end{lemma}

\begin{proof}
A test compactly supported in $U\cap\operatorname{Int}(M)$
extends smoothly by zero and is a global test.
This proves restriction. If the section vanishes on an open set,
the weak identity against tests in that set is zero; uniqueness
of functions representing distributions proves locality.

In a chart, the distributional Leibniz rule of
Lemma~\ref{lem:esquinas-cortes-extension-cero} gives the product
rule for the first order connection. Applying it successively
to the tensor bundles gives a sum over the choices of slots
differentiating $\chi$ and those differentiating $\mathbf u$.
For every $t$, there are $\binom{s}{t}$ choices.
Each term is a permutation of
$\nabla^{s-t}\chi\otimes\nabla_w^t\mathbf u$ preserving the order
of the slots in each factor. Product connections commute with
these permutations; induction proves the formula by separating,
at the next step, differentiation of the first and second factors.
For the induced tensor metrics, permutations are isometries,
and the norm of a tensor product is the product of the norms.
The pointwise inequality follows from the triangle inequality.
Taking $L^p$ norms and summing derivative orders, or taking maxima
when $p=\infty$, proves the multiplier estimate.

For zero extension, choose $\eta\in C_c^\infty(U)$ equal to one
near $K$. Locality makes the derivatives of $\mathbf v$ vanish
outside $K$ as well. For every interior test $\boldsymbol\sigma$,
\[
 \begin{aligned}
 \int_M\langle Z_U(\nabla_w^s\mathbf v),\boldsymbol\sigma\rangle
                                    \,d\lambda_{\mathbf g}
 &=\int_U\langle\nabla_w^s\mathbf v,\eta\boldsymbol\sigma\rangle
                                    \,d\lambda_{\mathbf g}\\
 &=\int_U\langle\mathbf v,(\nabla^s)_h^*(\eta\boldsymbol\sigma)\rangle
                                    \,d\lambda_{\mathbf g}\\
 &=\int_M\langle Z_U\mathbf v,(\nabla^s)_h^*\boldsymbol\sigma\rangle
                                    \,d\lambda_{\mathbf g}.
 \end{aligned}
\]
The product $\eta\boldsymbol\sigma$ is a local test, and near
$K$ its adjoint image agrees with the adjoint image of
$\boldsymbol\sigma$. These justify the two equalities.
The weak identity and preservation of the $L^p$ norms prove
the last assertion.
\end{proof}

\begin{proposition}[Completeness on arbitrary manifolds with corners]
\label{prop:esquinas-sobolev-completitud-global}
For every manifold with corners, compact or not, and every
$m\in\mathbb N_0$, $1\leq p\leq\infty$,
$W^{m,p}(M,\mathbf E)$ is a Banach space.
The map
\[
 J_m\mathbf u:=(\nabla_w^s\mathbf u)_{s=0}^m
\]
is an isometry into the product of the $L^p$ spaces, with the
power-sum norm for $p<\infty$ and the maximum norm for
$p=\infty$, and its image is closed.
\end{proposition}

\begin{proof}
The $L^p$ spaces of each tensor bundle are complete here as well.
To see this, choose a countable cover by orthonormal frames and
the disjoint measurable partition
$B_\alpha:=U_\alpha\setminus\bigcup_{\beta<\alpha}U_\beta$.
On $B_\alpha$, use the frame from $U_\alpha$.
The resulting components are measurable and isometrically
identify sections with measurable functions valued in a finite
dimensional Euclidean or Hermitian space. Their $L^p$ space is
complete by scalar completeness and equivalence of finite norms.

Suppose that $J_m\mathbf u_k\to(\mathbf v_s)_{s=0}^m$.
For an interior test $\boldsymbol\sigma$ and an order $s$,
Hölder's inequality bounds the differences of the integrals by
\[
 \|\nabla_w^s\mathbf u_k-\mathbf v_s\|_{L^p}
                      \|\boldsymbol\sigma\|_{L^{p'}},
 \qquad
 \|\mathbf u_k-\mathbf v_0\|_{L^p}
                      \|(\nabla^s)_h^*\boldsymbol\sigma\|_{L^{p'}},
 \qquad \frac{1}{p}+\frac{1}{p'}=1.
\]
The norms are over $M$ in the corresponding bundles.
The test and its adjoint image have compact interior support,
so their norms are finite even at the endpoints $p'=1$ and
$p'=\infty$. Passing to the limit in the weak identity gives
\[
 \int_M\langle\mathbf v_s,\boldsymbol\sigma\rangle\,d\lambda_{\mathbf g}
 =\int_M\langle\mathbf v_0,(\nabla^s)_h^*\boldsymbol\sigma\rangle
                                                   \,d\lambda_{\mathbf g}.
\]
Thus $\nabla_w^s\mathbf v_0=\mathbf v_s$ for every order,
and the image of $J_m$ is closed. It is an isometry by the
definition of the norm, and a closed subspace of a finite product
of Banach spaces is complete.
\end{proof}

\begin{definition}[Completion of sections smooth up to the faces]
\label{def:completacion-sobolev-haces-esquinas}
For $m\in\mathbb N_0$ and $1\leq p<\infty$, put
\[
 C^{m,p}(M,\mathbf E):=\Gamma(M,\mathbf E)\cap W^{m,p}(M,\mathbf E),
 \qquad
 H^{m,p}(M,\mathbf E):=
      \overline{C^{m,p}(M,\mathbf E)}^{\,\|\cdot\|_{W^{m,p}(M,\mathbf E)}}.
\]
The bar denotes completion, isometrically identified with its
closed image in the weak space.
\end{definition}

\begin{theorem}[Meyers--Serrin on manifolds with corners]
\label{teo:meyers-serrin-haces-esquinas}
Let $M$ be a Riemannian manifold with corners, with the data
on $\mathbf E$ smooth up to the faces. For every
$m\in\mathbb N_0$ and $1\leq p<\infty$,
\[
 H^{m,p}(M,\mathbf E)=W^{m,p}(M,\mathbf E)
\]
isometrically. Compactness, metric completeness, and uniform
geometric bounds are not required.
\end{theorem}

\begin{proof}
The completeness proposition gives inclusion of the completion
into the weak space. Fix $\mathbf u$ in the latter.
Proposition~\ref{prop:esquinas-cartas-regulares-cubierta} gives a
countable locally finite cover, frames extending to larger domains,
and cutoffs $\psi_\alpha$ compactly supported in each $U_\alpha$.
The model may have any depth $r_\alpha$. Put
$\Omega_\alpha:=\phi_\alpha(U_\alpha)$,
$\Omega_{\alpha,+}:=\Omega_\alpha\cap O_{n,r_\alpha}$, and
$h:=\operatorname{rank}\mathbf E$.
Write
\[
 \psi_\alpha\mathbf u=\sum_{a=1}^h u_\alpha^a\mathbf e_{\alpha,a},
 \quad f_\alpha^a:=u_\alpha^a\circ\phi_\alpha^{-1},
 \quad K_\alpha:=\phi_\alpha(\operatorname{supp}\psi_\alpha).
\]
The multiplier and local characterization lemmas show that
$f_\alpha^a\in W^{m,p}(\Omega_{\alpha,+})$, with essential support
in the same compact set $K_\alpha\Subset\Omega_\alpha$.
Lemma~\ref{lem:esquinas-aproximacion-ortante-soporte} supplies
smooth Euclidean approximants $f_{\alpha,k}^a$ converging to these
components and supported, after restriction to the model, in a
fixed compact set $L_\alpha\subseteq\Omega_\alpha$.
Define
\[
 \mathbf u_{\alpha,k}
    :=Z_{U_\alpha}\left(
       \sum_{a=1}^h(f_{\alpha,k}^a\circ\phi_\alpha)\mathbf e_{\alpha,a}
                    \right).
\]
These are globally smooth sections: their supports stay away from
the relative edge of the chart, and every component is a restriction
of a smooth Euclidean function.
Norm equivalence and extension of sections give
\[
 \|\mathbf u_{\alpha,k}-\psi_\alpha\mathbf u\|_{W^{m,p}(M,\mathbf E)}
 \leq C_\alpha\sum_{a=1}^h
       \|f_{\alpha,k}^a-f_\alpha^a\|_{W^{m,p}(\Omega_{\alpha,+})}
 \longrightarrow0.
\]
Here and in this norm, the approximants are understood as restricted
to the open domain.

For $q\in\mathbb N$, choose $k_{\alpha,q}$ so that the preceding
error is smaller than $2^{-(\alpha+q)}$, and put
$\mathbf u_q:=\sum_\alpha\mathbf u_{\alpha,k_{\alpha,q}}$.
The sum is locally finite and smooth. If
$\mathbf e_{\alpha,q}:=\psi_\alpha\mathbf u-
                                \mathbf u_{\alpha,k_{\alpha,q}}$,
the partial sums of errors satisfy
\[
 \left\|\sum_{\alpha=L+1}^N\mathbf e_{\alpha,q}
                         \right\|_{W^{m,p}(M,\mathbf E)}
 \leq\sum_{\alpha=L+1}^N2^{-(\alpha+q)}.
\]
They are Cauchy and converge in the Banach space.
On each neighborhood meeting only finitely many cover domains,
large partial sums agree with $\mathbf u-\mathbf u_q$.
Convergence in $L^p$ and a countable cover by such neighborhoods
identify the global limit with this difference.
Consequently,
\[
 \|\mathbf u-\mathbf u_q\|_{W^{m,p}(M,\mathbf E)}
 \leq\sum_{\alpha=1}^{\infty}2^{-(\alpha+q)}=2^{-q},
 \qquad
 \|\mathbf u_q\|_{W^{m,p}(M,\mathbf E)}
       \leq\|\mathbf u\|_{W^{m,p}(M,\mathbf E)}+2^{-q}.
\]
Thus every smooth approximant has finite norm and converges to
$\mathbf u$. For a finite cover, retain only the indices present;
the same bounds remain valid. This proves the equality.
\end{proof}

\begin{proposition}[Localization, completeness, and density]
\label{prop:sobolev-esquinas-localizacion-densidad}
Let $M$ be compact and put $h:=\operatorname{rank}\mathbf E$.
There are finitely many linear localization
maps $T_i$, with components supported in fixed compact sets
$K_i\subseteq G_{n,r_i}$, and reconstruction maps $S_i$ such that,
for every $m\in\mathbb N_0$ and $1\leq p\leq\infty$,
\[
 \mathbf u=\sum_{i=1}^s S_iT_i\mathbf u,\qquad
 C^{-1}\|\mathbf u\|_{W^{m,p}(M,\mathbf E)}
 \leq\sum_{i=1}^s\|T_i\mathbf u\|_{W^{m,p}(O_{n,r_i};\mathbb K^h)}
 \leq C\|\mathbf u\|_{W^{m,p}(M,\mathbf E)}.
\]
The maps $S_i$ are continuous at the same orders and exponents.
Global norms obtained from different smooth positive metrics
and connections are equivalent. The space is Banach, and sections
smooth up to the faces are dense for $p<\infty$.
\end{proposition}

\begin{proof}
Choose a finite cover of regular charts, frames, and cutoffs
$\psi_i$ as in Chapter 1. Write
$\Omega_i:=\phi_i(U_i)$. In each chart, put
\[
 \psi_i\mathbf u=\sum_{a=1}^h u_i^a\mathbf e_{i,a},\qquad
 T_i\mathbf u:=\bigl(Z_{\Omega_i}(u_i^a\circ\phi_i^{-1})\bigr)_{a=1}^h,
 \qquad K_i:=\phi_i(\operatorname{supp}\psi_i).
\]
The local characterization and scalar zero extension lemma identify
all their derivatives. Multiplication by $\psi_i$ and these
estimates give
\[
 \|T_i\mathbf u\|_{W^{m,p}(O_{n,r_i};\mathbb K^h)}
 \leq A_i\|\psi_i\mathbf u\|_{W^{m,p}(M,\mathbf E)}
 \leq A_iB_i\|\mathbf u\|_{W^{m,p}(M,\mathbf E)}.
\]
Choose $\xi_i\in C_c^\infty(\mathbb R^n)$ equal to one near $K_i$,
with $\operatorname{supp}\xi_i\cap G_{n,r_i}\Subset\Omega_i$,
and define
\[
 S_iv:=Z_{U_i}\left(
       \sum_{a=1}^h((\xi_iv^a)\circ\phi_i)\mathbf e_{i,a}\right).
\]
Scalar multiplication, norm equivalence on the chart, and
covariant extension prove
\[
 \|S_iv\|_{W^{m,p}(M,\mathbf E)}
        \leq D_i\|v\|_{W^{m,p}(O_{n,r_i};\mathbb K^h)}.
\]
Since $\xi_i=1$ near $K_i$, we have
$S_iT_i\mathbf u=\psi_i\mathbf u$. Summing gives reconstruction
and the lower bound; summing the upper estimates proves the
other inequality.

The maps $T_i$ use only charts, frames, and cutoffs.
For any other choice of metrics and connections, the local
characterization compares the covariant norm with these same
ordinary norms. The two comparisons prove both equality of
the spaces and equivalence of their norms.
Completeness has already been proved on arbitrary manifolds.
Meyers--Serrin gives density; here every smooth section has finite
norm because $M$ is compact, the derivatives are continuous up
to the faces, and the total volume is finite.
\end{proof}

\begin{remark}[Interior supports and boundary conditions]
\label{obs:sobolev-esquinas-W0}
For $m\in\mathbb N_0$ and $1\leq p\leq\infty$, define
$W_0^{m,p}(M,\mathbf E)$ as the closure, in $W^{m,p}(M,\mathbf E)$,
of the sections in $\Gamma_c(\operatorname{Int}(M),\mathbf E)$
extended by zero to $M$. Each such extension is smooth:
its support is a compact subset of the interior.
Theorem~\ref{teo:meyers-serrin-haces-esquinas} approximates, for
$p<\infty$, by smooth sections whose supports may reach the faces;
it does not assert that $W_0^{m,p}=W^{m,p}$.

The difference can be verified without a trace theorem.
On $M=[0,1]$, the constant function $1$ belongs to $W^{1,p}(M)$
for every $1\leq p\leq\infty$. If there were
$\varphi_\nu\in C_c^\infty((0,1))$ with
$\varphi_\nu\to1$ in $W^{1,p}((0,1))$, integration by parts would give
\[
 \int_0^1 t\varphi_\nu'(t)\,dt
 =-\int_0^1\varphi_\nu(t)\,dt.
\]
Let $p'$ be the conjugate exponent, with $p'=\infty$ if $p=1$ and
$p'=1$ if $p=\infty$. Hölder's inequality gives
\[
 \left|\int_0^1 t\varphi_\nu'(t)\,dt\right|
 \leq \|t\|_{L^{p'}(0,1)}\|\varphi_\nu'\|_{L^p(0,1)}
 \longrightarrow0,
 \qquad
 \left|\int_0^1(\varphi_\nu(t)-1)\,dt\right|
 \leq\|\varphi_\nu-1\|_{L^p(0,1)}
 \longrightarrow0.
\]
The preceding identity would require $0=-1$, a contradiction.
Thus $1\notin W_0^{1,p}([0,1])$, including the endpoints $p=1$
and $p=\infty$.
\end{remark}

\begin{theorem}[Sobolev and Rellich--Kondrashov on manifolds with corners]
\label{teo:sobolev-rellich-haces-esquinas}
Let $M$ be a compact Riemannian manifold with corners of dimension
$n\geq1$, and let $\mathbf E\to M$ be a finite rank bundle with data
smooth up to the faces. Let $\ell>j\geq0$ be integers and
$1\leq p,q<\infty$. If
\[
 \ell-\frac{n}{p}\geq j-\frac{n}{q},
\]
the embedding
$W^{\ell,p}(M,\mathbf E)\hookrightarrow W^{j,q}(M,\mathbf E)$
is continuous. If the inequality is strict, it is compact.
The same assertions hold for
$W_0^{\ell,p}(M,\mathbf E)\hookrightarrow W_0^{j,q}(M,\mathbf E)$.
\end{theorem}

\begin{proof}
Fix the maps $T_i,S_i$ from
Proposition~\ref{prop:sobolev-esquinas-localizacion-densidad},
the models $O_{n,r_i}$, and the compact sets $K_i$ supporting
the components of $T_i\mathbf u$.
By Lemma~\ref{lem:esquinas-encajes-localizados-ortante}, applied
with $\Omega=G_{n,r_i}$, every component belongs to
$W^{j,q}(O_{n,r_i})$, and the finite estimates give
\[
 \|T_i\mathbf u\|_{W^{j,q}(O_{n,r_i};\mathbb K^h)}
 \leq A_i\|T_i\mathbf u\|_{W^{\ell,p}(O_{n,r_i};\mathbb K^h)}.
\]
Continuity of $S_i$ at order $j$ and exponent $q$ places every
$S_iT_i\mathbf u$ in the global target space.
The reconstruction identity shows that their sum is $\mathbf u$.
Consequently,
\[
 \begin{aligned}
 \|\mathbf u\|_{W^{j,q}(M,\mathbf E)}
 &\leq\sum_{i=1}^sD_i
           \|T_i\mathbf u\|_{W^{j,q}(O_{n,r_i};\mathbb K^h)}\\
 &\leq\sum_{i=1}^sD_iA_i
           \|T_i\mathbf u\|_{W^{\ell,p}(O_{n,r_i};\mathbb K^h)}\\
 &\leq C\|\mathbf u\|_{W^{\ell,p}(M,\mathbf E)}.
 \end{aligned}
\]

For compactness, let $(\mathbf u_k)$ be bounded in the source space.
The localization estimates bound all components of $T_i\mathbf u_k$,
and their supports remain in $K_i$.
The compactness part of the scalar lemma gives a subsequence
converging in $W^{j,q}(O_{n,r_i})$ for each pair $(i,a)$.
There are finitely many pairs; successive extractions retain a
single subsequence for all of them. Let $v_i^a$ be their limits
and $v_i=(v_i^a)_{a=1}^h$.
Define $\mathbf u_\infty:=\sum_iS_iv_i$. Then
\[
 \|\mathbf u_k-\mathbf u_\infty\|_{W^{j,q}(M,\mathbf E)}
 \leq\sum_{i=1}^sD_i
      \|T_i\mathbf u_k-v_i\|_{W^{j,q}(O_{n,r_i};\mathbb K^h)}
 \longrightarrow0.
\]
This proves global compactness.

If $\mathbf u\in W_0^{\ell,p}$, choose
$\mathbf u_b\in\Gamma_c(\operatorname{Int}(M),\mathbf E)$
converging to $\mathbf u$ in the source norm.
The continuous estimate just proved gives convergence in $W^{j,q}$,
and the limit belongs to $W_0^{j,q}$ by the definition of closure.
Thus the inclusion restricts to the spaces with subscript zero.
Every bounded sequence in the source subspace has the compact
subsequence already constructed; closedness of the target
subspace preserves the limit. This proves its compactness as well.
\end{proof}

\begin{definition}[Hölder regularity up to corners]
\label{def:holder-haces-esquinas}
For $k\in\mathbb N_0$ and $0<\alpha<1$, use on the model the norm
\[
 \|f\|_{C^{k,\alpha}(G_{n,r})}
 :=\max_{|\beta|\leq k}\sup_{x\in G_{n,r}}|D^\beta f(x)|
  +\max_{|\beta|=k}\sup_{\substack{x,y\in G_{n,r}\\x\ne y}}
       \frac{|D^\beta f(x)-D^\beta f(y)|}{|x-y|^\alpha}.
\]
Write $[f]_{k,\alpha}$ for the second term. Derivatives on the
faces are derivatives of a local $C^k$ extension.
For the localization charts and cutoffs, define
\[
 \|\mathbf u\|_{k,\alpha,\mathcal A}
 :=\sum_{i=1}^s\sum_{a=1}^h
       \|(T_i\mathbf u)^a\|_{C^{k,\alpha}(G_{n,r_i})}.
\]
The space $\Gamma^{k,\alpha}(M,\mathbf E)$ consists of sections
of class $C^k$ up to the faces for which this norm is finite.
For $\alpha=0$, use only the first term of the model norm,
obtaining $\Gamma^k(M,\mathbf E)$.
\end{definition}

\begin{lemma}[Finite extension, independence, and Hölder completeness]
\label{lem:holder-esquinas-independencia}
Functions of class $C^k$ on $O_{n,r}$ whose derivatives through
order $k$ extend continuously to $G_{n,r}$ have local Euclidean
$C^k$ extensions. If their $C^{k,\alpha}(G_{n,r})$ norm is finite,
there is a bounded linear extension into
$C^{k,\alpha}(\mathbb R^n)$.
The norms in the preceding definition are equivalent for different
families of charts, frames, and cutoffs.
The spaces $\Gamma^{k,\alpha}$ and $\Gamma^k$ are Banach.
\end{lemma}

\begin{proof}
To extend one normal coordinate, solve
\[
 \sum_{\ell=1}^{k+1}a_\ell(-\ell)^b=1,\qquad 0\leq b\leq k,
\]
and use the reflection formula from
Lemma~\ref{lem:extension-sobolev-ortante-esquinas}.
Here derivatives are matched through order $k$.
We verify that the continuous extensions of these derivatives
are the derivatives of the glued function.
For a normal coordinate, the difference quotient at the face
is the average of the derivative on the normal segment, by
the fundamental theorem of calculus; continuity of that derivative
makes the average tend to its value at the face from both sides.
For a tangential coordinate, apply the same theorem on segments
with positive normal coordinate and let the latter tend to zero.
Uniform convergence on the compact set containing the segment
permits passage of the limit through the integral.
This gives the derivative identity on the face as well.
Applying these arguments to every derivative of order less than
$k$ proves by induction that the glued function is $C^k$ and
has the stated jets. The remaining normal coordinates may stay
nonnegative throughout; repeating the argument removes one face
at a time. For local data, first multiply by a cutoff compactly
supported in the relative domain.

We make the Hölder estimate explicit.
On either side, an order $k$ derivative is a finite combination
of derivatives of $f$ composed with dilated reflections.
The chain rule gives, for each reflection $R_\ell$,
\[
 [D^\beta f\circ R_\ell]_\alpha
             \leq \ell^\alpha[D^\beta f]_\alpha.
\]
Additional differentiation factors are finite powers of $\ell$
and are absorbed into a constant.
If $x$ and $y$ lie on opposite sides of the face, their segment
crosses it at a point $z$. The matching already proved gives
\[
 |D^\beta Ef(x)-D^\beta Ef(y)|
 \leq C[f]_{k,\alpha}(|x-z|^\alpha+|y-z|^\alpha)
 \leq2C[f]_{k,\alpha}|x-y|^\alpha.
\]
The segments remain in the respective convex models.
Uniform norms of the other derivatives are controlled directly
by the finite formula. Iteration yields
\[
 \mathcal R_{k,r}\colon C^{k,\alpha}(G_{n,r})
          \longrightarrow C^{k,\alpha}(\mathbb R^n),\qquad
 \|\mathcal R_{k,r}f\|_{C^{k,\alpha}(\mathbb R^n)}
          \leq C\|f\|_{C^{k,\alpha}(G_{n,r})}.
\]
The same construction without a seminorm gives the $C^k$ extension.

To compare localizations, use reconstruction
$\mathbf u=\sum_iS_iT_i\mathbf u$ and a second family $T_j'$.
On a chart overlap, the components have the form
\[
 (T_j'S_iv)^b(y)=\sum_{a=1}^h c_{ji}^{ba}(y)\,
                        v^a(F_{ji}(y)),\qquad
 F_{ji}=\phi_i\circ(\phi_j')^{-1}.
\]
The coefficients $c_{ji}^{ba}$ are products of cutoffs and frame
transition matrices. Their support is a compact subset of the
overlap, so the expression extends by zero before meeting any
artificial edge. Cover this compact set by finitely many balls
where the transition has a smooth Euclidean extension with bounded
derivatives. Use $\mathcal R_{k,r_i}v$ to compute compositions.
The chain and Leibniz rules express every derivative of order at
most $k$ as a finite sum of bounded smooth coefficients times
derivatives of $v$ of order at most $k$.
For the seminorms, use
\[
 [ab]_\alpha\leq\|a\|_\infty[b]_\alpha+\|b\|_\infty[a]_\alpha,
 \qquad
 [a\circ F]_\alpha\leq(\operatorname{Lip}F)^\alpha[a]_\alpha.
\]
Derivatives of order less than $k$ are Lipschitz by the bound on
their next derivatives and the fundamental theorem of calculus.
The smooth extensions of transitions and coefficients are Lipschitz
on the chosen balls. A finite subcover controls all sufficiently
close pairs. For the remaining pairs, at distance at least
$\delta>0$, the quotient is bounded by
$2\delta^{-\alpha}$ times the uniform norm.
For every pair of charts, this gives
\[
 \|T_j'S_iv\|_{C^{k,\alpha}(G_{n,r_j'};\mathbb K^h)}
       \leq C_{ji}\|v\|_{C^{k,\alpha}(G_{n,r_i};\mathbb K^h)}.
\]
Summing over pairs proves one norm comparison; interchanging
the families proves the other. For $\alpha=0$, omit only the
seminorm estimates.

For completeness, a Cauchy sequence on the model has uniform
limits of all derivatives through order $k$.
The fundamental theorem of calculus on interior segments identifies
these limits as derivatives of the limit function.
Continuity extends them to the faces, and the finite extension
just constructed gives $C^k$ regularity there.
If $[f_b-f_c]_{k,\alpha}<\varepsilon$ for large $b,c$, pass to the
limit in $c$ in each quotient and then take the supremum to get
$[f_b-f]_{k,\alpha}\leq\varepsilon$.
Thus convergence includes the seminorm.
Limits of localized components preserve their fixed compact
support. For a sequence $(\mathbf u_b)$ in the global space,
let $v_i$ be these limits of $T_i\mathbf u_b$, and set
$\mathbf u:=\sum_iS_iv_i$.
The preceding composition estimates imply
\[
 \|\mathbf u_b-\mathbf u\|_{k,\alpha,\mathcal A}
 \leq C\sum_i\|T_i\mathbf u_b-v_i\|_{C^{k,\alpha}(G_{n,r_i};\mathbb K^h)}
 \longrightarrow0.
\]
The reconstructed section is $C^k$ up to the faces.
For $\alpha=0$, the argument using uniform derivative norms
proves completeness of $\Gamma^k$.
\end{proof}

\begin{theorem}[Sobolev--Morrey up to every face]
\label{teo:sobolev-morrey-haces-esquinas}
Let $M$ be compact, with the preceding geometric data.
Let $m\in\mathbb N$, $1\leq p<\infty$, $k\in\mathbb N_0$,
and $0<\alpha<1$. If
\[
 k+\alpha\leq m-\frac{n}{p},
\]
every element of $W^{m,p}(M,\mathbf E)$ has a unique
representative in $\Gamma^{k,\alpha}(M,\mathbf E)$, and
\[
 \|\mathbf u\|_{k,\alpha,\mathcal A}
                  \leq C\|\mathbf u\|_{W^{m,p}(M,\mathbf E)}.
\]
The embedding is compact under the strict inequality.
If $m-\displaystyle\frac{n}{p}>k$, the embedding into
$\Gamma^k(M,\mathbf E)$ is continuous and compact.
\end{theorem}

\begin{proof}
Localize with $T_i$ and extend by
$V_i\mathbf u:=\mathcal E_{m,r_i}T_i\mathbf u$.
The components are supported in a fixed compact set
$\widehat K_i\Subset B_i$, where $B_i$ is a bounded open
Euclidean ball. Mollification with radius smaller than the
distance to $\partial B_i$ shows that
$V_i\mathbf u|_{B_i}\in W_0^{m,p}(B_i;\mathbb K^h)$,
as in the localized embedding lemma. Moreover,
\[
 \sum_{i=1}^s\|V_i\mathbf u\|_{W^{m,p}(\mathbb R^n;\mathbb K^h)}
                    \leq C\|\mathbf u\|_{W^{m,p}(M,\mathbf E)}.
\]
Theorem~\ref{teo:sobolev-morrey-euclidiano} provides representatives
$v_i^a$ of class $C^{k,\alpha}$ on $\overline B_i$, with the
corresponding bounds.

We specify their global extension.
Since the Sobolev class vanishes in $B_i\setminus\widehat K_i$,
every continuous derivative of the representative vanishes
there: it represents a weak derivative that is zero on that
open set. In particular, all derivatives through order $k$
vanish near $\partial B_i$. Its zero extension is therefore
$C^k$ on $\mathbb R^n$.
For the Hölder quotient, if both points lie in $\overline B_i$,
use the norm there. If one is outside and the derivative at the
other is nonzero, their distance is at least
$\operatorname{dist}(\widehat K_i,\partial B_i)>0$.
The uniform norm controls that quotient. We obtain
\[
 \|v_i^a\|_{C^{k,\alpha}(\mathbb R^n)}
           \leq C_i\|(V_i\mathbf u)^a\|_{W^{m,p}(B_i)}.
\]
The restrictions of these representatives to $G_{n,r_i}$
are reconstructed by $S_i$.
The section $\mathbf u^\sharp:=\sum_iS_i(v_i|_{G_{n,r_i}})$
agrees almost everywhere with $\mathbf u$.
The localization and reconstruction composition estimates of
the preceding lemma give
\[
 \|\mathbf u^\sharp\|_{k,\alpha,\mathcal A}
 \leq C\sum_{i=1}^s\sum_{a=1}^h
                  \|v_i^a\|_{C^{k,\alpha}(\mathbb R^n)}
 \leq C'\|\mathbf u\|_{W^{m,p}(M,\mathbf E)}.
\]

If two continuous representatives agree almost everywhere,
their difference vanishes on $\operatorname{Int}(M)$:
a nonzero value there would give an open set of positive
measure where it is nonzero, by continuity and positivity of
the volume density. The interior is dense and the difference
is continuous, so it vanishes on all faces as well.
This proves uniqueness.

Under the strict inequality, for a bounded sequence use the
same operators and balls.
Theorem~\ref{teo:rellich-kondrashov-holder-euclidiano} gives
subsequences converging in $C^{k,\alpha}(\overline B_i)$ for
every component. Finite extraction gives a single subsequence
for all of them. The limits and differences preserve support
in $\widehat K_i$, so the zero extension estimate just proved
preserves convergence in $C^{k,\alpha}(\mathbb R^n)$.
Continuous reconstruction and the finite sum give convergence
in $\Gamma^{k,\alpha}(M,\mathbf E)$.
Finally, if $m-\displaystyle\frac{n}{p}>k$, choose
$0<\alpha<1$ with $k+\alpha<m-\displaystyle\frac{n}{p}$.
The compact embedding already proved and the continuous inclusion
$\Gamma^{k,\alpha}\hookrightarrow\Gamma^k$ prove the last assertion.
\end{proof}

\begin{corollary}[Embeddings with infinite exponents at corners]
\label{cor:esquinas-encajes-exponentes-infinitos}
Under the preceding compact hypotheses, for $\ell>j\geq0$:
if $p<\infty$ and $\ell-\displaystyle\frac{n}{p}>j$, the embedding
$W^{\ell,p}(M,\mathbf E)\hookrightarrow W^{j,\infty}(M,\mathbf E)$
is compact. For $p=1$ and $\ell-j=n$, it is continuous.
For $p=\infty$ and every $1\leq q\leq\infty$, the embedding
$W^{\ell,\infty}(M,\mathbf E)\hookrightarrow W^{j,q}(M,\mathbf E)$
is compact.
\end{corollary}

\begin{proof}
In the first case, Morrey gives the compact embedding into
$\Gamma^j$. The coordinate expressions of covariant derivatives,
with coefficients bounded on the finite cover, prove continuity
of $\Gamma^j\hookrightarrow W^{j,\infty}$.
Their composition is compact.
At the endpoint $p=1$, $\ell-j=n$, localize and extend the
components to $W_0^{\ell,1}(B_i)$, with the supports and bounds
from the preceding proofs.
Corollary~\ref{cor:encajes-exponentes-infinitos-euclidianos}
gives the continuous estimate in $W^{j,\infty}(B_i)$.
Zero extension preserves derivatives because the supports stay
away from $\partial B_i$. Restriction and the maps $S_i$,
continuous for $p=\infty$, give the global bound.

If the source exponent is infinite, choose
$1<a<\infty$ with $\ell-\displaystyle\frac{n}{a}>j$.
The volume of $M$ is finite, and
\[
 \|\mathbf u\|_{W^{\ell,a}(M,\mathbf E)}
 \leq(\ell+1)^{\frac{1}{a}}
       \lambda_{\mathbf g}(M)^{\frac{1}{a}}
       \|\mathbf u\|_{W^{\ell,\infty}(M,\mathbf E)}.
\]
This inclusion is continuous, and the first case of the corollary
gives compactness from $W^{\ell,a}$ into $W^{j,\infty}$.
For $q<\infty$, the same volume estimate, at order $j$, gives
a continuous inclusion of $W^{j,\infty}$ into $W^{j,q}$.
The compositions prove all cases.
\end{proof}

\begin{example}[Embeddings and face compatibility on the square]
\label{ej:sobolev-cuadrado-compatibilidad-caras}
On $Q=[0,1]^2$ with the Euclidean metric, for every $2<p<\infty$ we have
\[
 W^{1,p}(Q)\hookrightarrow C^{0,\alpha}(Q),
 \qquad 0<\alpha\leq1-\frac{2}{p},
\]
with compactness if $\alpha<1-\displaystyle\frac{2}{p}$.
Also, $W^{2,2}(Q)\hookrightarrow C^{0,\alpha}(Q)$ for every
$0<\alpha<1$, and all these embeddings are compact. Values of the
continuous representative are defined even at the four vertices.

Let $u$ be one of these representatives, and set
$g_1(t):=u(t,0)$, $g_2(t):=u(0,t)$. They are continuous and necessarily
\[
 g_1(0)=u(0,0)=g_2(0).
\]
Thus prescribing $g_1\equiv0$ and $g_2\equiv1$ on two intersecting
sides cannot arise from an element of these Sobolev spaces. Although
both data are smooth on each side separately, their values are
incompatible at the vertex. The assertion uses the stated continuous
regularity: for orders and exponents without that embedding, a
Sobolev class is not assumed to have a point value at a vertex.
\end{example}

\begin{corollary}[Localized embeddings without compactness of the manifold]
\label{cor:esquinas-encajes-locales-cortes}
Let $M$ be a Riemannian manifold with corners, possibly noncompact,
and let $\mathbf E\to M$ be a finite-rank bundle with the smooth data
of Definition~\ref{def:sobolev-haces-esquinas}. Let
$\chi,\eta\in C_c^\infty(M)$, with $\eta=1$ on a neighborhood of
$\operatorname{supp}\chi$, and let $\ell>j\geq0$ and $1\leq p,q<\infty$.
If $\ell-\displaystyle\frac{n}{p}\geq j-\displaystyle\frac{n}{q}$,
there exists $C_{\chi,\eta}>0$ such that
\[
 \|\chi\mathbf u\|_{W^{j,q}(M,\mathbf E)}
 \leq C_{\chi,\eta}\|\eta\mathbf u\|_{W^{\ell,p}(M,\mathbf E)}
 \qquad
 \bigl(\mathbf u\in W^{\ell,p}_{\operatorname{loc}}(M,\mathbf E)\bigr).
\]
If the index inequality is strict and
$\{\|\eta\mathbf u_\nu\|_{W^{\ell,p}(M,\mathbf E)}\}_\nu$ is bounded,
then $\{\chi\mathbf u_\nu\}_\nu$ has a subsequence converging
in $W^{j,q}(M,\mathbf E)$.
\end{corollary}

\begin{proof}
The case $\chi=0$ is immediate; assume $K:=\operatorname{supp}\chi$
is nonempty. By Lemma~\ref{lem:esquinas-localidad-cortes-secciones},
$\eta\mathbf u$ belongs to the global space $W^{\ell,p}$:
in each chart meeting its support, its components and derivatives
are obtained through the weak Leibniz rule; a finite cover of the
compact support and local norm equivalence give global
integrability. Outside that support all its weak derivatives
vanish by locality.

Cover $K$ by regular charts $(U_i,\phi_i)$, $1\leq i\leq N$,
with frames $(\mathbf e_{i,a})_{a=1}^h$ defined on a neighborhood
of $\overline U_i$, as in
Proposition~\ref{prop:esquinas-cartas-regulares-cubierta}.
A partition of unity subordinate to
$U_1,\dots,U_N,M\setminus K$ provides functions
$\psi_i\in C_c^\infty(U_i)$ such that
$\displaystyle\sum_{i=1}^N\psi_i=1$ on $K$.
Put $\Omega_i:=\phi_i(U_i)\subset G_{n,r_i}$ and
$K_i:=\phi_i(\operatorname{supp}(\psi_i\chi))$.
These are compact subsets of $\Omega_i$.
If $K_i$ is empty, the corresponding contribution is zero and
can be removed from the family. Write
\[
 \psi_i\chi\eta\mathbf u
 =\displaystyle\sum_{a=1}^h f_i^a\mathbf e_{i,a},
 \qquad
 v_i^a:=Z_{\Omega_i}(f_i^a\circ\phi_i^{-1})
 \in W^{\ell,p}(O_{n,r_i}).
\]
The Leibniz rule, Lemma~\ref{lem:esquinas-componentes-normas-locales},
and Lemma~\ref{lem:esquinas-cortes-extension-cero} give constants
$A_i>0$, independent of $\mathbf u$, such that
\[
 \displaystyle\sum_{a=1}^h\|v_i^a\|_{W^{\ell,p}(O_{n,r_i})}
 \leq A_i\|\eta\mathbf u\|_{W^{\ell,p}(M,\mathbf E)}.
\]
All supports lie in the same $K_i$.
Lemma~\ref{lem:esquinas-encajes-localizados-ortante} gives
\[
 \displaystyle\sum_{a=1}^h\|v_i^a\|_{W^{j,q}(O_{n,r_i})}
 \leq B_i\displaystyle\sum_{a=1}^h
             \|v_i^a\|_{W^{\ell,p}(O_{n,r_i})}.
\]

To reconstruct the sections, choose
$\xi_i\in C_c^\infty(\mathbb R^n)$ equal to $1$ near $K_i$,
with $\operatorname{supp}\xi_i\cap G_{n,r_i}$ compact in $\Omega_i$.
For $v=(v^1,\dots,v^h)$ put
\[
 S_i v:=
 Z_{U_i}\left(
 \displaystyle\sum_{a=1}^h
       ((\xi_i v^a)\circ\phi_i)\mathbf e_{i,a}
 \right).
\]
The product is interpreted on $\Omega_i\cap O_{n,r_i}$.
Its support lies in the fixed compact set
$\operatorname{supp}\xi_i\cap G_{n,r_i}\subset\Omega_i$.
The Leibniz rule and local norm equivalence show that
the section in parentheses belongs to $W^{j,q}(U_i,\mathbf E)$.
Lemma~\ref{lem:esquinas-localidad-cortes-secciones} permits
its extension by zero to $M$ and gives, for some $D_i>0$,
\[
 \|S_i v\|_{W^{j,q}(M,\mathbf E)}
 \leq D_i\displaystyle\sum_{a=1}^h
                     \|v^a\|_{W^{j,q}(O_{n,r_i})}.
\]
Since $\eta=1$ near $K$, we have
$\chi\mathbf u=\displaystyle\sum_{i=1}^N S_i v_i$.
The preceding three estimates and the triangle inequality
prove the embedding with
$C_{\chi,\eta}:=\displaystyle\sum_{i=1}^N D_i B_i A_i$.

In the strict case, each component sequence $v_{i,\nu}^a$
is bounded in $W^{\ell,p}(O_{n,r_i})$ and supported in the
fixed compact set $K_i$. By the compact part of the localized lemma,
successively extract subsequences for the $Nh$ pairs $(i,a)$.
The last subsequence satisfies
$v_{i,\nu}^a\to v_i^a$ in $W^{j,q}(O_{n,r_i})$ for all of them.
Then $\mathbf w:=\displaystyle\sum_{i=1}^N S_i v_i$
belongs to $W^{j,q}(M,\mathbf E)$ and
\[
 \|\chi\mathbf u_\nu-\mathbf w\|_{W^{j,q}(M,\mathbf E)}
 \leq\displaystyle\sum_{i=1}^N D_i
        \displaystyle\sum_{a=1}^h
        \|v_{i,\nu}^a-v_i^a\|_{W^{j,q}(O_{n,r_i})}
 \longrightarrow0.
\]
Only constants from a finite family of charts covering $K$
have been used; no geometric uniformity on all of $M$ was assumed.
\end{proof}

\begin{remark}[Local scope, traces, and boundary problems]
\label{obs:alcance-analisis-esquinas}
Corollary~\ref{cor:esquinas-encajes-locales-cortes} proves localized embeddings and compactness for noncompact manifolds. Its constants depend on the cutoffs and charts covering their supports. A uniform global bound requires control of those geometric data; the definition of a manifold with corners alone does not provide that control.

The embedding and compactness proofs do not require the total
boundary to be a smooth manifold. The trace theorems on $\partial M$
for smooth boundary do use that structure, and must be reformulated
in face charts with the corresponding compatibility conditions.
Where only one face is active, the reduced chart is a half-space
and the established local trace theorem applies. Near an intersection,
extension to Euclidean space permits restrictions to individual
faces, but describing their joint image requires additional conditions,
as the preceding example illustrates. The embeddings proved here do
not, by themselves, assert elliptic regularity up to a corner for
Dirichlet or Neumann problems, or automatically enlarge the bounded
geometry results in later chapters.
\end{remark}

\section{The weak Kato inequality}

Earlier we proved Kato's inequality for smooth sections. We now extend it to sections in $W^{m,p}(M,\mathbf{E})$ by Theorem~\ref{teo:meyers-serrin-haz-frontera}. When $1<p<\infty$, reflexivity of the $L^p(M,\mathbf{F})$ spaces of sections of finite-rank bundles allows us to extract the weakly convergent subsequences needed. To include the case $p=1$, where $L^1(M,\mathbf{F})$ is not reflexive, we use the Dunford--Pettis theorem.

Recall that, if $(X,\mathcal A,\mu)$ is a measure space, a family $\mathcal F\subseteq L^1(X,\mathcal A,\mu)$ is uniformly integrable if $\displaystyle\sup_{f\in\mathcal F}\|f\|_{L^1(X)}<\infty$ and, for every $\varepsilon>0$, there is $\delta>0$ such that $\displaystyle\int_A|f|\,d\mu<\varepsilon$ for any $f\in\mathcal F$ and $A\in\mathcal A$ with $\mu(A)<\delta$.

\begin{theorem}[Dunford--Pettis]\label{dunford pettis}\index{Dunford Pettis@Dunford--Pettis} Let $(X,\mathcal A,\mu)$ be a finite measure space. A family $\mathcal F\subseteq L^1(X,\mathcal A,\mu)$ has weakly compact closure in $L^1(X,\mathcal A,\mu)$ if and only if $\mathcal F$ is uniformly integrable. \end{theorem}

This is the classical form of the Dunford--Pettis criterion for a finite measure; see \cite{Bogachev2007}.

By Theorem~\ref{teo:dualidad-reflexividad-Lp}, the weak topology of $L^1(X,\mathcal A,\mu)$ agrees with $\sigma(L^1(X,\mathcal A,\mu),L^\infty(X,\mathcal A,\mu))$. Thus the preceding statement is equivalent to the classical formulation using $L^1(X,\mathcal A,\mu)$--$L^\infty(X,\mathcal A,\mu)$ duality.

\begin{theorem}[Dunford--Pettis for vector bundles] \label{dunford pettis para haces vectoriales}\index{Dunford Pettis for vector bundles@Dunford--Pettis for vector bundles} Let $(M,\mathbf{g})$ be a compact Riemannian manifold with or without boundary, and let $\mathbf{E}\to M$ be a smooth real vector bundle of finite rank equipped with a bundle metric $\mathbf{h}_{\mathbf{E}}$. For a family $\mathcal F\subseteq L^1(M,\mathbf{E})$, the following assertions are equivalent: \begin{enumerate}[label=(\alph*)] \item The weak closure of $\mathcal F$ in $L^1(M,\mathbf{E})$ is compact.

\item The family $|\mathcal F|_{\mathbf{h}_{\mathbf{E}}}:=\{|\mathbf{u}|_{\mathbf{h}_{\mathbf{E}}}\mid \mathbf{u}\in\mathcal F\}
 \subseteq L^1(M)$ is uniformly integrable. \end{enumerate} \end{theorem}

\begin{proof} Since $M$ is compact, there is a finite open cover $(U_\alpha)_{\alpha=1}^N$ of $M$ with a local trivialization of $\mathbf{E}$ defined over each $U_\alpha$, together with a smooth partition of unity $(\chi_\alpha)_{\alpha=1}^N$ subordinate to this cover. On each $U_\alpha$ fix a smooth local orthonormal frame $(\mathbf{e}_{\alpha,1},\dots,\mathbf{e}_{\alpha,r})$.

For $\alpha\in\{1,\dots,N\}$ and $a\in\{1,\dots,r\}$, define \[
T_{\alpha,a}\colon L^1(M,\mathbf{E})\longrightarrow L^1(U_\alpha),
\qquad
T_{\alpha,a}(\mathbf{u})
:=
\langle\chi_\alpha \mathbf{u},\mathbf{e}_{\alpha,a}\rangle_{\mathbf{h}_{\mathbf{E}}}.
\] Since \[
|T_{\alpha,a}(\mathbf{u})|
\leq
\chi_\alpha|\mathbf{u}|_{\mathbf{h}_{\mathbf{E}}}
\leq
|\mathbf{u}|_{\mathbf{h}_{\mathbf{E}}},
\] we have \[
\|T_{\alpha,a}(\mathbf{u})\|_{L^1(U_\alpha)}
\leq
\|\mathbf{u}\|_{L^1(M,\mathbf{E})}.
\] Thus $T_{\alpha,a}$ is linear and norm continuous. By Proposition~\ref{prop:operador-lineal-continuo-topologias-debiles}, it is also continuous for the weak topologies.

First suppose that the weak closure $C:=\overline{\mathcal F}^{\,\sigma(L^1(M,\mathbf{E}),(L^1(M,\mathbf{E}))')}$ is compact. Weak continuity of $T_{\alpha,a}$ implies that $T_{\alpha,a}(C)$ is weakly compact in $L^1(U_\alpha)$ and contains $T_{\alpha,a}(\mathcal F)$. By the Hahn--Banach theorem, $(L^1(U_\alpha))'$ separates points of $L^1(U_\alpha)$; by Proposition~\ref{prop: topologia debil dualidad seminormas}, the weak topology of this space is Hausdorff. Consequently, $T_{\alpha,a}(C)$ is weakly closed. Thus \[
\overline{T_{\alpha,a}(\mathcal F)}^{\,\sigma(L^1(U_\alpha),(L^1(U_\alpha))')}
\subseteq
T_{\alpha,a}(C),
\] and the weak closure of $T_{\alpha,a}(\mathcal F)$ is compact. Theorem~\ref{dunford pettis} shows that each family $T_{\alpha,a}(\mathcal F)$ is uniformly integrable.

For every $\mathbf{u}\in\mathcal F$, since \[
\mathbf{u}=
\sum_{\alpha=1}^N\chi_\alpha \mathbf{u}
\] and, on $U_\alpha$, \[
\chi_\alpha \mathbf{u}
=
\sum_{a=1}^rT_{\alpha,a}(\mathbf{u})\mathbf{e}_{\alpha,a},
\] we have \[
|\mathbf{u}|_{\mathbf{h}_{\mathbf{E}}}
\leq
\sum_{\alpha=1}^N|\chi_\alpha \mathbf{u}|_{\mathbf{h}_{\mathbf{E}}}
\leq
\sum_{\alpha=1}^N\sum_{a=1}^r|T_{\alpha,a}(\mathbf{u})|.
\] In particular, \[
\sup_{\mathbf{u}\in\mathcal F}\|\mathbf{u}\|_{L^1(M,\mathbf{E})}
\leq
\sum_{\alpha=1}^N\sum_{a=1}^r
\sup_{\mathbf{u}\in\mathcal F}
\|T_{\alpha,a}(\mathbf{u})\|_{L^1(U_\alpha)}
<\infty.
\] Given $\varepsilon>0$, for each $(\alpha,a)\in\{1,\dots,N\}\times\{1,\dots,r\}$ there is $\delta_{\alpha,a}>0$ such that \[
\int_A|T_{\alpha,a}(\mathbf{u})|\,d\lambda_{\mathbf{g}}
<
\frac{\varepsilon}{Nr}
\] for every $\mathbf{u}\in\mathcal F$ and every measurable set $A\subseteq U_\alpha$ with $\lambda_{\mathbf{g}}(A)<\delta_{\alpha,a}$. Taking \[
\delta
:=
\min_{\substack{1\leq\alpha\leq N\\1\leq a\leq r}}
\delta_{\alpha,a},
\] if $A\subseteq M$ is measurable and $\lambda_{\mathbf{g}}(A)<\delta$, then \[
\int_A|\mathbf{u}|_{\mathbf{h}_{\mathbf{E}}}\,d\lambda_{\mathbf{g}}
\leq
\sum_{\alpha=1}^N\sum_{a=1}^r
\int_{A\cap U_\alpha}|T_{\alpha,a}(\mathbf{u})|\,d\lambda_{\mathbf{g}}
<
\varepsilon
\] for every $\mathbf{u}\in\mathcal F$. Consequently, $|\mathcal F|_{\mathbf{h}_{\mathbf{E}}}$ is uniformly integrable.

Conversely, suppose that $|\mathcal F|_{\mathbf{h}_{\mathbf{E}}}$ is uniformly integrable. Since \[
|T_{\alpha,a}(\mathbf{u})|\leq|\mathbf{u}|_{\mathbf{h}_{\mathbf{E}}},
\] we have \[
\sup_{\mathbf{u}\in\mathcal F}
\|T_{\alpha,a}(\mathbf{u})\|_{L^1(U_\alpha)}
\leq
\sup_{\mathbf{u}\in\mathcal F}\|\mathbf{u}\|_{L^1(M,\mathbf{E})}
<\infty.
\] Moreover, given $\varepsilon>0$, there is $\delta>0$ such that, if $A\subseteq U_\alpha$ is measurable and $\lambda_{\mathbf{g}}(A)<\delta$, then \[
\int_A|T_{\alpha,a}(\mathbf{u})|\,d\lambda_{\mathbf{g}}
\leq
\int_A|\mathbf{u}|_{\mathbf{h}_{\mathbf{E}}}\,d\lambda_{\mathbf{g}}
<
\varepsilon
\] for every $\mathbf{u}\in\mathcal F$. Thus each family $T_{\alpha,a}(\mathcal F)$ is uniformly integrable. Theorem~\ref{dunford pettis} implies that its weak closure \[
K_{\alpha,a}
:=
\overline{T_{\alpha,a}(\mathcal F)}^{\,\sigma(L^1(U_\alpha),(L^1(U_\alpha))')}
\] is compact.

Consider the finite product \[
\mathbb X
:=
\prod_{\alpha=1}^N\prod_{a=1}^rL^1(U_\alpha)
\] equipped, for example, with the norm \[
\|(f_{\alpha,a})_{\alpha,a}\|_{\mathbb X}
:=
\sum_{\alpha=1}^N\sum_{a=1}^r
\|f_{\alpha,a}\|_{L^1(U_\alpha)}.
\] Since the product has finitely many factors, the topology of this norm agrees with the product topology. By Proposition~\ref{prop:topologia-debil-producto-finito}, the weak topology of $\mathbb X$ agrees with the product of the weak topologies of the factors. Consequently, \[
K
:=
\prod_{\alpha=1}^N\prod_{a=1}^rK_{\alpha,a}
\] is weakly compact in $\mathbb X$.

Define \[
T\colon L^1(M,\mathbf{E})\longrightarrow\mathbb X,
\qquad
T(\mathbf{u}):=(T_{\alpha,a}(\mathbf{u}))_{\alpha,a},
\] and \[
S\colon \mathbb X\longrightarrow L^1(M,\mathbf{E}),
\qquad
S((f_{\alpha,a})_{\alpha,a})
:=
\sum_{\alpha=1}^N\sum_{a=1}^r
\widetilde{f_{\alpha,a}\mathbf{e}_{\alpha,a}},
\] where the tilde denotes extension by zero outside $U_\alpha$. The operator $S$ is linear and satisfies \[
\|S((f_{\alpha,a})_{\alpha,a})\|_{L^1(M,\mathbf{E})}
\leq
\sum_{\alpha=1}^N\sum_{a=1}^r
\|f_{\alpha,a}\|_{L^1(U_\alpha)}
=
\|(f_{\alpha,a})_{\alpha,a}\|_{\mathbb X}.
\] Thus $S$ is norm continuous and, by Proposition~\ref{prop:operador-lineal-continuo-topologias-debiles}, is continuous from the weak topology of $\mathbb X$ to the weak topology of $L^1(M,\mathbf{E})$. Since $K$ is weakly compact, $S(K)$ is weakly compact in $L^1(M,\mathbf{E})$.

Finally, for every $\mathbf{u}\in L^1(M,\mathbf{E})$, \[
S(T(\mathbf{u}))=
\sum_{\alpha=1}^N\sum_{a=1}^r
\widetilde{
\langle\chi_\alpha \mathbf{u},\mathbf{e}_{\alpha,a}\rangle_{\mathbf{h}_{\mathbf{E}}}
\mathbf{e}_{\alpha,a}}=
\sum_{\alpha=1}^N\chi_\alpha \mathbf{u}
=
\mathbf{u}
\] almost everywhere. Since $T(\mathcal F)\subseteq K$, it follows that $\mathcal F=S(T(\mathcal F))\subseteq S(K)$. Again by Hahn--Banach and Proposition~\ref{prop: topologia debil dualidad seminormas}, the weak topology of $L^1(M,\mathbf{E})$ is Hausdorff. Thus the compact set $S(K)$ is weakly closed. Consequently, the weak closure of $\mathcal F$ is contained in $S(K)$ and is a weakly closed subset of this compact set. Therefore, the weak closure of $\mathcal F$ is compact. \end{proof}

We need the following consequence of weak lower semicontinuity of convex functionals.

\begin{lemma}\label{semicontinuidad norma ponderada} Let $(M,\mathbf{g})$ be a Riemannian manifold with or without boundary, let $\mathbf{F}\to M$ be a vector bundle of finite rank equipped with a bundle metric $\mathbf{h}_{\mathbf{F}}$, let $1\leq p<\infty$, and let $p'$ be the conjugate exponent of $p$. For every nonnegative function $\varphi\in L^{p'}(M)$, the functional \[
\Phi_\varphi\colon L^p(M,\mathbf{F})\longrightarrow[0,\infty),
\qquad
\Phi_\varphi(\mathbf{Z})
:=
\int_M\varphi|\mathbf{Z}|_{\mathbf{h}_{\mathbf{F}}}\,d\lambda_{\mathbf{g}},
\] is convex, norm continuous, and lower semicontinuous with respect to $\sigma(L^p(M,\mathbf{F}),(L^p(M,\mathbf{F}))')$. \end{lemma}

\begin{proof} Convexity follows from the triangle inequality in the fibers. Moreover, for any $\mathbf{Z},\mathbf{W}\in L^p(M,\mathbf{F})$, Hölder's inequality (Proposition~\ref{desigualdad de holder}) implies \[
|\Phi_\varphi(\mathbf{Z})-\Phi_\varphi(\mathbf{W})|
\leq
\int_M\varphi|\mathbf{Z}-\mathbf{W}|_{\mathbf{h}_{\mathbf{F}}}\,d\lambda_{\mathbf{g}}\leq
\|\varphi\|_{L^{p'}(M)}
\|\mathbf{Z}-\mathbf{W}\|_{L^p(M,\mathbf{F})}.
\] Thus $\Phi_\varphi$ is norm continuous. Since it is also convex, Corollary~\ref{cor: funcional convexo semicontinuo es debilmente semicontinuo} implies that it is lower semicontinuous with respect to $\sigma(L^p(M,\mathbf{F}),(L^p(M,\mathbf{F}))')$. \end{proof} By Lemma \ref{lema para kato}, we know that, for a smooth section $\mathbf{u}$, $|\mathbf{u}|_{\mathbf{h}_{\mathbf{E}}}$ admits a weak derivative almost everywhere on $M$. Now that we have the concept of a \textit{weak covariant derivative}, it is natural to ask whether this almost-everywhere defined classical derivative agrees with a weak derivative. The following proposition answers this question, which will be important in proving our weak version of Kato's inequality. \begin{lemma}[Classical and weak derivatives of the modulus of a smooth section] \label{lem:derivada-debil-modulo-seccion-suave} Let $(M,\mathbf{g})$ be a Riemannian manifold with or without boundary, and let $\mathbf{F}\to M$ be a smooth real vector bundle of finite rank equipped with a bundle metric $\mathbf{h}_{\mathbf{F}}$ and a compatible connection $\nabla^{\mathbf{F}}$. Let $1\leq p<\infty$ and $\mathbf{v}\in\Gamma(\mathbf{F})$. Then $|\mathbf{v}|_{\mathbf{h}_{\mathbf{F}}}\in W_{\mathrm{loc}}^{1,p}(M)$ and $\nabla_w|\mathbf{v}|_{\mathbf{h}_{\mathbf{F}}}
=
\nabla|\mathbf{v}|_{\mathbf{h}_{\mathbf{F}}}$ almost everywhere on $M$, where the right-hand side denotes the classical covariant derivative defined almost everywhere. Moreover, $\left|\nabla_w|\mathbf{v}|_{\mathbf{h}_{\mathbf{F}}}\right|_{\mathbf{g}}
\leq
|\nabla^{\mathbf{F}}\mathbf{v}|_{\mathbf{g},\mathbf{h}_{\mathbf{F}}}$ almost everywhere. If, in addition, $\mathbf{v}\in L^p(M,\mathbf{F})$ and $\nabla^{\mathbf{F}}\mathbf{v}\in L^p(M,T^*M\otimes \mathbf{F})$, then $|\mathbf{v}|_{\mathbf{h}_{\mathbf{F}}}\in W^{1,p}(M)$. In particular, this last conclusion holds automatically when $M$ is compact. \end{lemma}

\begin{proof} By Lemma~\ref{lema para kato}, $\nabla|\mathbf{v}|_{\mathbf{h}_{\mathbf{F}}}$ exists almost everywhere on $M$. Define the measurable $1$-form $G$ by \[
G(x):=
\begin{cases}
\nabla|\mathbf{v}|_{\mathbf{h}_{\mathbf{F}}}(x),
&
|\mathbf{v}(x)|_{\mathbf{h}_{\mathbf{F}}}>0,
\\[1ex]
0,
&
|\mathbf{v}(x)|_{\mathbf{h}_{\mathbf{F}}}=0.
\end{cases}
\] On the set where $|\mathbf{v}|_{\mathbf{h}_{\mathbf{F}}}>0$, the classical Kato inequality, Theorem~\ref{desigualdad de kato}, applied with $s=0$, gives $|G|_{\mathbf{g}}
=
\left|\nabla|\mathbf{v}|_{\mathbf{h}_{\mathbf{F}}}\right|_{\mathbf{g}}
\leq
|\nabla^{\mathbf{F}}\mathbf{v}|_{\mathbf{g},\mathbf{h}_{\mathbf{F}}}$. On the set where $|\mathbf{v}|_{\mathbf{h}_{\mathbf{F}}}=0$, this inequality also holds by the definition of $G$. Consequently, $|G|_{\mathbf{g}}\leq|\nabla^{\mathbf{F}}\mathbf{v}|_{\mathbf{g},\mathbf{h}_{\mathbf{F}}}$ almost everywhere on $M$.

For each $\varepsilon>0$, consider the smooth function \[
f_\varepsilon
:=
\sqrt{|\mathbf{v}|_{\mathbf{h}_{\mathbf{F}}}^{\,2}+\varepsilon^2}-\varepsilon.
\] On the open set where $|\mathbf{v}|_{\mathbf{h}_{\mathbf{F}}}>0$, the chain rule implies \[
\nabla f_\varepsilon
=
\frac{|\mathbf{v}|_{\mathbf{h}_{\mathbf{F}}}}
{\sqrt{|\mathbf{v}|_{\mathbf{h}_{\mathbf{F}}}^{\,2}+\varepsilon^2}}
\nabla|\mathbf{v}|_{\mathbf{h}_{\mathbf{F}}}
=
\frac{|\mathbf{v}|_{\mathbf{h}_{\mathbf{F}}}}
{\sqrt{|\mathbf{v}|_{\mathbf{h}_{\mathbf{F}}}^{\,2}+\varepsilon^2}}G.
\] If $|\mathbf{v}|_{\mathbf{h}_{\mathbf{F}}}=0$ at a point, the smooth function $|\mathbf{v}|_{\mathbf{h}_{\mathbf{F}}}^{\,2}$ attains its minimum there, and hence $\nabla|\mathbf{v}|_{\mathbf{h}_{\mathbf{F}}}^{\,2}=0$. Consequently, \[
\nabla f_\varepsilon
=
\frac{\nabla|\mathbf{v}|_{\mathbf{h}_{\mathbf{F}}}^{\,2}}
{2\sqrt{|\mathbf{v}|_{\mathbf{h}_{\mathbf{F}}}^{\,2}+\varepsilon^2}}
=
0
\] at that point. Thus throughout $M$ we have \[
\nabla f_\varepsilon
=
\frac{|\mathbf{v}|_{\mathbf{h}_{\mathbf{F}}}}
{\sqrt{|\mathbf{v}|_{\mathbf{h}_{\mathbf{F}}}^{\,2}+\varepsilon^2}}G.
\]

Since \[
0\leq
\frac{|\mathbf{v}|_{\mathbf{h}_{\mathbf{F}}}}
{\sqrt{|\mathbf{v}|_{\mathbf{h}_{\mathbf{F}}}^{\,2}+\varepsilon^2}}
\leq1,
\] it follows that $|\nabla f_\varepsilon|_{\mathbf{g}}
\leq|G|_{\mathbf{g}}
\leq|\nabla^{\mathbf{F}}\mathbf{v}|_{\mathbf{g},\mathbf{h}_{\mathbf{F}}}$.

Moreover, $\nabla f_\varepsilon\longrightarrow G$ as $\varepsilon\to 0^{+}$ pointwise on $M$, since on the set where $|\mathbf{v}|_{\mathbf{h}_{\mathbf{F}}}>0$ the preceding quotient converges to $1$, while on the set where $|\mathbf{v}|_{\mathbf{h}_{\mathbf{F}}}=0$ both $1$-forms vanish. Let $K\subseteq M$ be compact. Since $\mathbf{v}$ is smooth, $|\nabla^{\mathbf{F}}\mathbf{v}|_{\mathbf{g},\mathbf{h}_{\mathbf{F}}}\in L^p(K)$ and, by the dominated convergence theorem, $\nabla f_\varepsilon\longrightarrow G$ in $L^p(K,T^*M)$.

On the other hand, $0\leq f_\varepsilon\leq|\mathbf{v}|_{\mathbf{h}_{\mathbf{F}}}$ and $f_\varepsilon\longrightarrow|\mathbf{v}|_{\mathbf{h}_{\mathbf{F}}}$ pointwise. Again by dominated convergence, $f_\varepsilon\longrightarrow|\mathbf{v}|_{\mathbf{h}_{\mathbf{F}}}$ in $L^p(K)$.

Now let $\boldsymbol{\eta}\in\Gamma_{c,\operatorname{Int}(M)}(T^*M)$ and choose a compact set $K$ containing the support of $\boldsymbol{\eta}$. Since $f_\varepsilon$ is smooth, \[
\int_M
\langle\nabla f_\varepsilon,\boldsymbol{\eta}\rangle_{\mathbf{g}}\,d\lambda_{\mathbf{g}}
=
\int_M
f_\varepsilon\nabla^*\boldsymbol{\eta}\,d\lambda_{\mathbf{g}}.
\] The preceding strong convergences and Hölder's inequality allow us to pass to the limit and obtain \[
\int_M\langle G,\boldsymbol{\eta}\rangle_{\mathbf{g}}\,d\lambda_{\mathbf{g}}
=
\int_M|\mathbf{v}|_{\mathbf{h}_{\mathbf{F}}}\nabla^*\boldsymbol{\eta}\,d\lambda_{\mathbf{g}}.
\] By the definition of weak covariant derivatives and local Sobolev spaces, \[
|\mathbf{v}|_{\mathbf{h}_{\mathbf{F}}}\in W_{\mathrm{loc}}^{1,p}(M)
\qquad\text{and}\qquad
\nabla_w|\mathbf{v}|_{\mathbf{h}_{\mathbf{F}}}=G
\] almost everywhere.

It remains to identify $G$ with the classical derivative. On the set where $|\mathbf{v}|_{\mathbf{h}_{\mathbf{F}}}>0$, this equality holds by the definition of $G$ itself. If $|\mathbf{v}|_{\mathbf{h}_{\mathbf{F}}}=0$ at a point where the classical derivative exists, the nonnegative function $|\mathbf{v}|_{\mathbf{h}_{\mathbf{F}}}$ attains its minimum there, so $\nabla|\mathbf{v}|_{\mathbf{h}_{\mathbf{F}}}=0=G$.

Since the classical derivative exists almost everywhere, we conclude that $\nabla_w|\mathbf{v}|_{\mathbf{h}_{\mathbf{F}}}
=
\nabla|\mathbf{v}|_{\mathbf{h}_{\mathbf{F}}}$ almost everywhere on $M$. The estimate already obtained for $G$ gives \[
\left|\nabla_w|\mathbf{v}|_{\mathbf{h}_{\mathbf{F}}}\right|_{\mathbf{g}}
\leq
|\nabla^{\mathbf{F}}\mathbf{v}|_{\mathbf{g},\mathbf{h}_{\mathbf{F}}}
\] almost everywhere.

Finally, if $\mathbf{v}\in L^p(M,\mathbf{F})$ and $\nabla^{\mathbf{F}}\mathbf{v}\in L^p(M,T^*M\otimes \mathbf{F})$, the preceding weak identity therefore implies that $|\mathbf{v}|_{\mathbf{h}_{\mathbf{F}}}\in W^{1,p}(M)$. \end{proof} \begin{theorem}[Weak Kato inequality] \label{desigualdad debil de kato}\index{weak Kato inequality@weak Kato inequality} Let $(M,\mathbf{g})$ be a compact Riemannian manifold with or without boundary, and let $\mathbf{E}\to M$ be a smooth real vector bundle of finite rank equipped with a bundle metric $\mathbf{h}_{\mathbf{E}}$ and a compatible connection $\nabla^{\mathbf{E}}$. Let $m\in\mathbb N$, $1\leq p<\infty$, and $\mathbf{u}\in W^{m,p}(M,\mathbf{E})$. Then, for every $s\in\{0,\dots,m-1\}$, $|\nabla_w^s\mathbf{u}|_{\mathbf{g},\mathbf{h}_{\mathbf{E}}}\in W^{1,p}(M)$ and \[
\left|
\nabla_w|\nabla_w^s\mathbf{u}|_{\mathbf{g},\mathbf{h}_{\mathbf{E}}}
\right|_{\mathbf{g}}
\leq
|\nabla_w^{s+1}\mathbf{u}|_{\mathbf{g},\mathbf{h}_{\mathbf{E}}}
\] almost everywhere on $M$, where $\nabla_w^0\mathbf{u}:=\mathbf{u}$. \end{theorem}

\begin{proof} Fix $s\in\{0,\dots,m-1\}$. By Theorem~\ref{meyers-serrin-haz} if $\partial M=\varnothing$, or Theorem~\ref{teo:meyers-serrin-haz-frontera} if $\partial M\neq\varnothing$, there is a sequence $(\mathbf{u}_k)_{k\in\mathbb N}\subseteq\Gamma(\mathbf{E})\cap W^{m,p}(M,\mathbf{E})$ such that $\mathbf{u}_k\to \mathbf{u}$ in $W^{m,p}(M,\mathbf{E})$. The reverse triangle inequality in the fibers implies \[
\left|
|\nabla^s\mathbf{u}_k|_{\mathbf{g},\mathbf{h}_{\mathbf{E}}}
-
|\nabla_w^s\mathbf{u}|_{\mathbf{g},\mathbf{h}_{\mathbf{E}}}
\right|
\leq
|\nabla^s\mathbf{u}_k-\nabla_w^s\mathbf{u}|_{\mathbf{g},\mathbf{h}_{\mathbf{E}}},
\] so that $|\nabla^s\mathbf{u}_k|_{\mathbf{g},\mathbf{h}_{\mathbf{E}}}\to|\nabla_w^s\mathbf{u}|_{\mathbf{g},\mathbf{h}_{\mathbf{E}}}$ in $L^p(M)$.

By Lemma~\ref{lem:derivada-debil-modulo-seccion-suave}, each function $|\nabla^s\mathbf{u}_k|_{\mathbf{g},\mathbf{h}_{\mathbf{E}}}$ admits a first weak covariant derivative agreeing almost everywhere with its classical derivative. Kato's inequality for smooth sections, Theorem~\ref{desigualdad de kato}, gives \[
\left|
\nabla_w|\nabla^s\mathbf{u}_k|_{\mathbf{g},\mathbf{h}_{\mathbf{E}}}
\right|_{\mathbf{g}}
\leq
|\nabla^{s+1}\mathbf{u}_k|_{\mathbf{g},\mathbf{h}_{\mathbf{E}}}
\] almost everywhere. Since $|\nabla^s\mathbf{u}_k|_{\mathbf{g},\mathbf{h}_{\mathbf{E}}}\in L^p(M)$ and $|\nabla^{s+1}\mathbf{u}_k|_{\mathbf{g},\mathbf{h}_{\mathbf{E}}}\in L^p(M)$, we conclude that $|\nabla^s\mathbf{u}_k|_{\mathbf{g},\mathbf{h}_{\mathbf{E}}}\in W^{1,p}(M)$.

Moreover, again by the reverse triangle inequality, \[
\left|
|\nabla^{s+1}\mathbf{u}_k|_{\mathbf{g},\mathbf{h}_{\mathbf{E}}}
-
|\nabla_w^{s+1}\mathbf{u}|_{\mathbf{g},\mathbf{h}_{\mathbf{E}}}
\right|
\leq
|\nabla^{s+1}\mathbf{u}_k-\nabla_w^{s+1}\mathbf{u}|_{\mathbf{g},\mathbf{h}_{\mathbf{E}}},
\] so that $|\nabla^{s+1}\mathbf{u}_k|_{\mathbf{g},\mathbf{h}_{\mathbf{E}}}\to
|\nabla_w^{s+1}\mathbf{u}|_{\mathbf{g},\mathbf{h}_{\mathbf{E}}}$ in $L^p(M)$. In particular, the sequence $\bigl(\nabla_w|\nabla^s\mathbf{u}_k|_{\mathbf{g},\mathbf{h}_{\mathbf{E}}}\bigr)_{k\in\mathbb N}$ is bounded in $L^p(M,T^*M)$.

If $1<p<\infty$, reflexivity of $L^p(M,T^*M)$ allows us to extract a subsequence, which we continue to denote in the same way, and find $\mathbf{Z}\in L^p(M,T^*M)$ such that \[
\nabla_w|\nabla^s\mathbf{u}_k|_{\mathbf{g},\mathbf{h}_{\mathbf{E}}}
\rightharpoonup \mathbf{Z}
\quad\text{in }L^p(M,T^*M).
\]

Now suppose that $p=1$. Since $|\nabla^{s+1}\mathbf{u}_k|_{\mathbf{g},\mathbf{h}_{\mathbf{E}}}\to
|\nabla_w^{s+1}\mathbf{u}|_{\mathbf{g},\mathbf{h}_{\mathbf{E}}}$ in $L^1(M)$, the family \[
\left\{
|\nabla^{s+1}\mathbf{u}_k|_{\mathbf{g},\mathbf{h}_{\mathbf{E}}}
\;\middle|\;
k\in\mathbb N
\right\}
\] is uniformly integrable. Convergence in $L^1(M)$ implies that this family is bounded in $L^1(M)$. On the other hand, given $\varepsilon>0$, there is $N\in\mathbb N$ such that \[
\left\|
|\nabla^{s+1}\mathbf{u}_k|_{\mathbf{g},\mathbf{h}_{\mathbf{E}}}
-
|\nabla_w^{s+1}\mathbf{u}|_{\mathbf{g},\mathbf{h}_{\mathbf{E}}}
\right\|_{L^1(M)}
<
\frac{\varepsilon}{2}
\] for every $k\geq N$. Since every function in $L^1(M)$ has an integral absolutely continuous with respect to $\lambda_{\mathbf{g}}$, and only the functions $|\nabla_w^{s+1}\mathbf{u}|_{\mathbf{g},\mathbf{h}_{\mathbf{E}}},
|\nabla^{s+1}\mathbf{u}_1|_{\mathbf{g},\mathbf{h}_{\mathbf{E}}},\dots,
|\nabla^{s+1}\mathbf{u}_{N-1}|_{\mathbf{g},\mathbf{h}_{\mathbf{E}}}$ are involved, there is $\delta>0$ such that $\lambda_{\mathbf{g}}(A)<\delta$ implies \[
\int_A|\nabla_w^{s+1}\mathbf{u}|_{\mathbf{g},\mathbf{h}_{\mathbf{E}}}\,d\lambda_{\mathbf{g}}
<
\frac{\varepsilon}{2}
\] and \[
\int_A|\nabla^{s+1}\mathbf{u}_j|_{\mathbf{g},\mathbf{h}_{\mathbf{E}}}\,d\lambda_{\mathbf{g}}
<
\varepsilon
\] for every $j\in\{1,\dots,N-1\}$. For $k\geq N$ we obtain \[
\int_A|\nabla^{s+1}\mathbf{u}_k|_{\mathbf{g},\mathbf{h}_{\mathbf{E}}}\,d\lambda_{\mathbf{g}}
\leq
\left\|
|\nabla^{s+1}\mathbf{u}_k|_{\mathbf{g},\mathbf{h}_{\mathbf{E}}}
-
|\nabla_w^{s+1}\mathbf{u}|_{\mathbf{g},\mathbf{h}_{\mathbf{E}}}
\right\|_{L^1(M)}
+
\int_A|\nabla_w^{s+1}\mathbf{u}|_{\mathbf{g},\mathbf{h}_{\mathbf{E}}}\,d\lambda_{\mathbf{g}}
<
\varepsilon.
\] Thus the preceding family is uniformly integrable. Since \[
\left|
\nabla_w|\nabla^s\mathbf{u}_k|_{\mathbf{g},\mathbf{h}_{\mathbf{E}}}
\right|_{\mathbf{g}}
\leq
|\nabla^{s+1}\mathbf{u}_k|_{\mathbf{g},\mathbf{h}_{\mathbf{E}}},
\] the family \[
\left\{
\left|
\nabla_w|\nabla^s\mathbf{u}_k|_{\mathbf{g},\mathbf{h}_{\mathbf{E}}}
\right|_{\mathbf{g}}
\;\middle|\;
k\in\mathbb N
\right\}
\] is also uniformly integrable. Theorem~\ref{dunford pettis para haces vectoriales}, applied to the bundle $T^*M$, implies that the weak closure of $\{\nabla_w|\nabla^s\mathbf{u}_k|_{\mathbf{g},\mathbf{h}_{\mathbf{E}}}\mid k\in\mathbb N\}$ in $L^1(M,T^*M)$ is compact. By the Eberlein--Šmulian theorem (Theorem~\ref{teo: eberlein smulian}), there are a subsequence and $\mathbf{Z}\in L^1(M,T^*M)$ such that \[
\nabla_w|\nabla^s\mathbf{u}_k|_{\mathbf{g},\mathbf{h}_{\mathbf{E}}}
\rightharpoonup \mathbf{Z}
\quad\text{in }L^1(M,T^*M).
\] Thus in both cases, after extracting a subsequence if necessary, we may assume that \[
\nabla_w|\nabla^s\mathbf{u}_k|_{\mathbf{g},\mathbf{h}_{\mathbf{E}}}
\rightharpoonup \mathbf{Z}
\quad\text{in }L^p(M,T^*M)
\] for some $\mathbf{Z}\in L^p(M,T^*M)$.

Let $\boldsymbol{\eta}\in\Gamma_{c,\operatorname{Int}(M)}(T^*M)$. Since $\nabla_w|\nabla^s\mathbf{u}_k|_{\mathbf{g},\mathbf{h}_{\mathbf{E}}}$ is the weak covariant derivative of $|\nabla^s\mathbf{u}_k|_{\mathbf{g},\mathbf{h}_{\mathbf{E}}}$, we have \[
\int_M
\left\langle
\nabla_w|\nabla^s\mathbf{u}_k|_{\mathbf{g},\mathbf{h}_{\mathbf{E}}},
\boldsymbol{\eta}
\right\rangle_{\mathbf{g}}\,d\lambda_{\mathbf{g}}
=
\int_M
|\nabla^s\mathbf{u}_k|_{\mathbf{g},\mathbf{h}_{\mathbf{E}}}\nabla^*\boldsymbol{\eta}\,d\lambda_{\mathbf{g}}.
\] Weak convergence of the left-hand side and strong convergence $|\nabla^s\mathbf{u}_k|_{\mathbf{g},\mathbf{h}_{\mathbf{E}}}\to|\nabla_w^s\mathbf{u}|_{\mathbf{g},\mathbf{h}_{\mathbf{E}}}$ in $L^p(M)$ allow us to pass to the limit, since $\nabla^*\boldsymbol{\eta}\in L^{p'}(M)$, with the interpretation $p'=\infty$ when $p=1$, and by Hölder's inequality (Proposition~\ref{desigualdad de holder}) \[
\left|
\int_M
\left(
|\nabla^s\mathbf{u}_k|_{\mathbf{g},\mathbf{h}_{\mathbf{E}}}
-
|\nabla_w^s\mathbf{u}|_{\mathbf{g},\mathbf{h}_{\mathbf{E}}}
\right)\nabla^*\boldsymbol{\eta}\,d\lambda_{\mathbf{g}}
\right|
\leq
\left\|
|\nabla^s\mathbf{u}_k|_{\mathbf{g},\mathbf{h}_{\mathbf{E}}}
-
|\nabla_w^s\mathbf{u}|_{\mathbf{g},\mathbf{h}_{\mathbf{E}}}
\right\|_{L^p(M)}
\|\nabla^*\boldsymbol{\eta}\|_{L^{p'}(M)}
\longrightarrow0.
\] Consequently, \[
\int_M\langle \mathbf{Z},\boldsymbol{\eta}\rangle_{\mathbf{g}}\,d\lambda_{\mathbf{g}}
=
\int_M
|\nabla_w^s\mathbf{u}|_{\mathbf{g},\mathbf{h}_{\mathbf{E}}}\nabla^*\boldsymbol{\eta}\,d\lambda_{\mathbf{g}}.
\] By the definition of weak covariant derivatives, $|\nabla_w^s\mathbf{u}|_{\mathbf{g},\mathbf{h}_{\mathbf{E}}}\in W^{1,p}(M)$ and $\mathbf{Z}=
\nabla_w|\nabla_w^s\mathbf{u}|_{\mathbf{g},\mathbf{h}_{\mathbf{E}}}$.

Finally, let $\varphi\in C_c^\infty(\operatorname{Int}(M))$ be nonnegative. By Lemma~\ref{semicontinuidad norma ponderada}, \[
\int_M
\varphi
\left|
\nabla_w|\nabla_w^s\mathbf{u}|_{\mathbf{g},\mathbf{h}_{\mathbf{E}}}
\right|_{\mathbf{g}}\,d\lambda_{\mathbf{g}}
\leq
\liminf_{k\to\infty}
\int_M
\varphi
\left|
\nabla_w|\nabla^s\mathbf{u}_k|_{\mathbf{g},\mathbf{h}_{\mathbf{E}}}
\right|_{\mathbf{g}}\,d\lambda_{\mathbf{g}}.
\] Kato's inequality for smooth sections and convergence $|\nabla^{s+1}\mathbf{u}_k|_{\mathbf{g},\mathbf{h}_{\mathbf{E}}}\to
|\nabla_w^{s+1}\mathbf{u}|_{\mathbf{g},\mathbf{h}_{\mathbf{E}}}$ in $L^p(M)$ imply \[
\liminf_{k\to\infty}
\int_M
\varphi
\left|
\nabla_w|\nabla^s\mathbf{u}_k|_{\mathbf{g},\mathbf{h}_{\mathbf{E}}}
\right|_{\mathbf{g}}\,d\lambda_{\mathbf{g}}
\leq
\lim_{k\to\infty}
\int_M
\varphi|\nabla^{s+1}\mathbf{u}_k|_{\mathbf{g},\mathbf{h}_{\mathbf{E}}}\,d\lambda_{\mathbf{g}}
=
\int_M
\varphi|\nabla_w^{s+1}\mathbf{u}|_{\mathbf{g},\mathbf{h}_{\mathbf{E}}}\,d\lambda_{\mathbf{g}},
\] where the last convergence follows from Hölder's inequality in Proposition~\ref{desigualdad de holder}, since $\varphi\in L^{p'}(M)$. Thus \[
\int_M
\varphi
\left(
\left|
\nabla_w|\nabla_w^s\mathbf{u}|_{\mathbf{g},\mathbf{h}_{\mathbf{E}}}
\right|_{\mathbf{g}}
-
|\nabla_w^{s+1}\mathbf{u}|_{\mathbf{g},\mathbf{h}_{\mathbf{E}}}
\right)
\,d\lambda_{\mathbf{g}}
\leq0
\] for every nonnegative $\varphi\in C_c^\infty(\operatorname{Int}(M))$. Since the function in parentheses belongs to $L^1_{\mathrm{loc}}(\operatorname{Int}(M))$, we conclude that \[
\left|
\nabla_w|\nabla_w^s\mathbf{u}|_{\mathbf{g},\mathbf{h}_{\mathbf{E}}}
\right|_{\mathbf{g}}
\leq
|\nabla_w^{s+1}\mathbf{u}|_{\mathbf{g},\mathbf{h}_{\mathbf{E}}}
\] almost everywhere on $\operatorname{Int}(M)$. Since $\partial M$ has measure zero, the inequality also holds almost everywhere on $M$. \end{proof}

\begin{remark}\label{obs:kato-debil-restriccion-debe-uso-teorema-teo-meyers} The restriction $p<\infty$ comes from the use of Theorem~\ref{teo:meyers-serrin-haz-frontera}, since smooth sections are not, in general, dense in $W^{m,\infty}(M,\mathbf{E})$ in its norm. \end{remark}

\begin{corollary}[Weak Kato inequality]\label{desigualdad debil de kato general}\index{weak Kato inequality@weak Kato inequality} Let $(M,\mathbf{g})$ be a Riemannian manifold with or without boundary, and let $\mathbf{E}\to M$ be a smooth real vector bundle of finite rank equipped with a bundle metric $\mathbf{h}_{\mathbf{E}}$ and a compatible connection $\nabla^{\mathbf{E}}$. Let $m\in\mathbb N$, $1\leq p<\infty$, and $\mathbf{u}\in W^{m,p}(M,\mathbf{E})$. Then, for every $s\in\{0,\dots,m-1\}$, we have $|\nabla_w^s\mathbf{u}|_{\mathbf{g},\mathbf{h}_{\mathbf{E}}}\in W^{1,p}(M)$ and $\displaystyle\left|\nabla_w|\nabla_w^s\mathbf{u}|_{\mathbf{g},\mathbf{h}_{\mathbf{E}}}\right|_{\mathbf{g}}\leq|\nabla_w^{s+1}\mathbf{u}|_{\mathbf{g},\mathbf{h}_{\mathbf{E}}}$ almost everywhere on $M$. \end{corollary}

\begin{proof} Apply a smooth exhaustion to the manifold $\operatorname{Int}(M)$, which has no boundary: choose increasing open sets $(\Omega_j)_{j\in\mathbb N}$ with smooth boundary, $\overline{\Omega_j}\Subset\operatorname{Int}(M)$ and $\operatorname{Int}(M)=\displaystyle\bigcup_{j\in\mathbb N}\Omega_j$. If $M$ is already compact, the preceding theorem also gives the conclusion directly; the exhaustion allows all cases to be treated by the same interior argument. By Proposition~\ref{prop: restriccion de seccion de Sobolev}, for every $j\in\mathbb N$ and every $r\leq m$ we have $\mathbf{u}\restriction_{\Omega_j}\in W^{m,p}(\Omega_j,\mathbf{E}\restriction_{\Omega_j})$ and $\nabla_w^r(\mathbf{u}\restriction_{\Omega_j})=(\nabla_w^r\mathbf{u})\restriction_{\Omega_j}$ almost everywhere on $\Omega_j$. Identifying $\Omega_j$ with the interior of the compact manifold with boundary $\overline{\Omega_j}$, Theorem~\ref{desigualdad debil de kato} implies that $|\nabla_w^s\mathbf{u}|_{\mathbf{g},\mathbf{h}_{\mathbf{E}}}\restriction_{\Omega_j}\in W^{1,p}(\Omega_j)$ and there is $\mathbf{v}_j\in L^p(\Omega_j,T^*\Omega_j)$ such that $\mathbf{v}_j=\nabla_w|\nabla_w^s\mathbf{u}|_{\mathbf{g},\mathbf{h}_{\mathbf{E}}}$ on $\Omega_j$ and $\displaystyle |\mathbf{v}_j|_{\mathbf{g}}\leq|\nabla_w^{s+1}\mathbf{u}|_{\mathbf{g},\mathbf{h}_{\mathbf{E}}}$ almost everywhere on $\Omega_j$.

If $j\leq k$, then $|\nabla_w^s\mathbf{u}|_{\mathbf{g},\mathbf{h}_{\mathbf{E}}}\restriction_{\Omega_j}$ is the restriction of $|\nabla_w^s\mathbf{u}|_{\mathbf{g},\mathbf{h}_{\mathbf{E}}}\restriction_{\Omega_k}$ to $\Omega_j$, and by Proposition~\ref{prop: restriccion de seccion de Sobolev} the weak derivative of this restriction is $\mathbf{v}_k\restriction_{\Omega_j}$, while by definition it is also $\mathbf{v}_j$. Uniqueness of the weak covariant derivative gives $\mathbf{v}_k\restriction_{\Omega_j}=\mathbf{v}_j$ almost everywhere on $\Omega_j$. Thus there is a measurable section $\mathbf{v}$ of $T^*M$ such that $\mathbf{v}\restriction_{\Omega_j}=\mathbf{v}_j$ almost everywhere on $\Omega_j$ for every $j\in\mathbb N$. Since $\displaystyle |\mathbf{v}|_{\mathbf{g}}\leq|\nabla_w^{s+1}\mathbf{u}|_{\mathbf{g},\mathbf{h}_{\mathbf{E}}}$ almost everywhere on each $\Omega_j$ and $\operatorname{Int}(M)=\displaystyle\bigcup_{j\in\mathbb N}\Omega_j$, this inequality holds almost everywhere on $M$, since $\partial M$ has zero $n$-dimensional measure. Extend $\mathbf v$ by any value, for example zero, on $\partial M$, giving $\displaystyle\int_M|\mathbf{v}|_{\mathbf{g}}^p\,d\lambda_{\mathbf{g}}\leq\int_M|\nabla_w^{s+1}\mathbf{u}|_{\mathbf{g},\mathbf{h}_{\mathbf{E}}}^p\,d\lambda_{\mathbf{g}}<\infty$ and consequently $\mathbf{v}\in L^p(M,T^*M)$.

Finally, let $\mathbf{X}\in\Gamma_{c,\operatorname{Int}(M)}(T^*M)$. Since $\supp(\mathbf{X})$ is a compact subset of the interior and $(\Omega_j)_{j\in\mathbb N}$ is an increasing exhaustion of it, there is $j\in\mathbb N$ such that $\supp(\mathbf{X})\subseteq\Omega_j$. Since $\mathbf{v}\restriction_{\Omega_j}=\mathbf{v}_j$ and $\mathbf{v}_j$ is the weak derivative of $|\nabla_w^s\mathbf{u}|_{\mathbf{g},\mathbf{h}_{\mathbf{E}}}\restriction_{\Omega_j}$, we have \[
\int_M\langle \mathbf{v},\mathbf{X}\rangle_{\mathbf{g}}\,d\lambda_{\mathbf{g}}
=\int_{\Omega_j}\langle \mathbf{v}_j,\mathbf{X}\rangle_{\mathbf{g}}\,d\lambda_{\mathbf{g}}
=\int_{\Omega_j}|\nabla_w^s\mathbf{u}|_{\mathbf{g},\mathbf{h}_{\mathbf{E}}}\nabla^*\mathbf{X}\,d\lambda_{\mathbf{g}}
=\int_M|\nabla_w^s\mathbf{u}|_{\mathbf{g},\mathbf{h}_{\mathbf{E}}}\nabla^*\mathbf{X}\,d\lambda_{\mathbf{g}}.
\] The first and last equalities follow because $\supp(\mathbf{X})\subseteq\Omega_j$ and, by locality of $\nabla^*$, $\supp(\nabla^*\mathbf{X})\subseteq\supp(\mathbf{X})$. Thus $\mathbf{v}=\nabla_w|\nabla_w^s\mathbf{u}|_{\mathbf{g},\mathbf{h}_{\mathbf{E}}}$. Since $|\nabla_w^s\mathbf{u}|_{\mathbf{g},\mathbf{h}_{\mathbf{E}}}\in L^p(M)$ and $\mathbf{v}\in L^p(M,T^*M)$, we conclude that $|\nabla_w^s\mathbf{u}|_{\mathbf{g},\mathbf{h}_{\mathbf{E}}}\in W^{1,p}(M)$, and the inequality $\displaystyle\left|\nabla_w|\nabla_w^s\mathbf{u}|_{\mathbf{g},\mathbf{h}_{\mathbf{E}}}\right|_{\mathbf{g}}\leq|\nabla_w^{s+1}\mathbf{u}|_{\mathbf{g},\mathbf{h}_{\mathbf{E}}}$ holds almost everywhere on $M$. \end{proof}

\section{The Poincaré inequality in the compact case}

The Poincaré inequality controls the $L^p(M,\mathbf E)$ norm of a section, and hence its Sobolev norm, solely in terms of its first weak covariant derivative. More precisely, an estimate of the form $\displaystyle \|\mathbf{u}\|_{W^{1,p}(M,\mathbf{E})}\leq C\|\nabla_w^{\mathbf{E}}\mathbf{u}\|_{L^p(M,T^*M\otimes \mathbf{E})}$ is equivalent, by the definition of the Sobolev norm, to estimating $\|\mathbf{u}\|_{L^p(M,\mathbf{E})}$ by $\|\nabla_w^{\mathbf{E}}\mathbf{u}\|_{L^p(M,T^*M\otimes \mathbf{E})}$. Such an estimate, however, cannot hold throughout $W^{1,p}(M,\mathbf{E})$ when the connection admits nonzero parallel sections.

The trivial bundle immediately exhibits this obstruction. If $M$ is a compact Riemannian manifold and $\mathbf{E}=M\times\mathbb R$ is equipped with the trivial metric and connection, and $u\equiv c$ with $c\neq0$, then $u\in W^{1,p}(M)$ and $\nabla_wu=0$, while $\|u\|_{L^p(M)}>0$. In a general vector bundle, the same obstruction is given by the kernel of the weak covariant derivative, that is, the sections $\mathbf{u}$ satisfying $\nabla_w^{\mathbf{E}}\mathbf{u}=0$. Thus the Poincaré inequality does not control a section absolutely, but only modulo the space of parallel sections.

This observation allows a unified treatment of the different forms of the inequality. One possibility is to project each section onto the kernel of the connection and estimate the difference between the section and its parallel component. For the trivial bundle over a connected manifold, the kernel consists of constant functions and the corresponding projection is precisely the mean. Another possibility is to impose a homogeneous Dirichlet condition: if every connected component of the manifold has nonempty boundary, no nonzero parallel section can have zero trace. Both formulations remove the same obstruction by different mechanisms.

We begin by describing the space of parallel sections and a continuous projection onto it. For this we need a preliminary lemma about functions whose weak derivative is $0$. \begin{lemma}[Functions with zero weak covariant derivative] \label{lem:derivada-debil-nula-funcion-constante} \index{function with zero weak covariant derivative@function with zero weak covariant derivative} Let $(M,\mathbf{g})$ be a Riemannian manifold with or without boundary, let $1\leq p<\infty$, and let $u\in W^{1,p}(M)$. If \[
\nabla_wu=0
\] almost everywhere on $M$, then $u$ is constant almost everywhere on each connected component of $M$. In particular, if $M$ is connected, there is $c\in\mathbb R$ such that $u=c$ almost everywhere on $M$. \end{lemma}

\begin{proof} If $\dim M=0$, the assertion is immediate. Henceforth suppose that $\dim M=n\geq 1$ and fix a connected component $M_\alpha$ of $M$.

By Proposition~\ref{bolacoordenadaregular}, $M_\alpha$ admits a countable open cover by regular coordinate balls and regular coordinate half-balls; write this cover as \[
(U_i,\phi_i)_{i\in I}.
\] For each $i\in I$, define \[
\Omega_i:=\phi_i\bigl(U_i\cap \operatorname{Int}(M)\bigr).
\] Then $\Omega_i$ is an open ball or an open half-ball in $\mathbb R^n$, as appropriate; in particular, $\Omega_i$ is connected.

We now examine the situation locally. Since $u\in W^{1,p}(M)$ and $\nabla_wu=0$ almost everywhere on $M$, Proposition~\ref{prop: restriccion de seccion de Sobolev} implies that \[
u\restriction_{U_i}\in W^{1,p}(U_i)
\quad\text{and}\quad
\nabla_w(u\restriction_{U_i})
=
(\nabla_wu)\restriction_{U_i}
=
0
\quad\text{almost everywhere on }U_i.
\] If $M$ has boundary, the preceding application is justified by the identification described in Section~\ref{sec:sobolev-haces-frontera}. Indeed, $u\restriction_{\operatorname{Int}(M)}$ belongs to $W^{1,p}(\operatorname{Int}(M))$ and its weak covariant derivative agrees almost everywhere with $(\nabla_wu)\restriction_{\operatorname{Int}(M)}=0$. Apply the proposition to the manifold without boundary $\operatorname{Int}(M)$ and the open set $U_i\cap\operatorname{Int}(M)=\operatorname{Int}(U_i)$. Using that identification again, now on $U_i$, gives
$u\restriction_{U_i}\in W^{1,p}(U_i)$ and
$\nabla_w(u\restriction_{U_i})=0$ almost everywhere on $U_i$.
Indeed, the test sections have compact support in
$\operatorname{Int}(U_i)$, and $U_i\cap\partial M$ has zero
Riemannian measure.

Since we are dealing with scalar functions, the covariant derivative agrees with the differential. Thus, applying Lemmas~\ref{lema:meyers-serrin-haz-E} and \ref{lema:meyers-serrin-haz-frontera} to the trivial rank-one bundle, we find that the local representation \[
v_i:=u\circ\phi_i^{-1}
\] belongs to $W^{1,p}(\Omega_i)$ and satisfies \[
D_jv_i=0
\quad\text{almost everywhere on }\Omega_i,
\qquad j=1,\dots,n.
\] In other words, all first Euclidean weak derivatives of $v_i$ vanish almost everywhere.

By Lemma~\ref{lem:derivadas-debiles-nulas-euclidiano}, there is a constant $c_i\in\mathbb R$ such that \[
v_i=c_i
\quad\text{almost everywhere on }\Omega_i.
\] Returning to the manifold, this means that \[
u=c_i
\quad\text{almost everywhere on }U_i\cap \operatorname{Int}(M).
\]

We claim that in fact also \[
u=c_i
\quad\text{almost everywhere on }U_i.
\] Indeed, if $U_i$ is a regular coordinate ball, then $U_i\subseteq \operatorname{Int}(M)$ and there is nothing to prove. If $U_i$ is a regular coordinate half-ball, then \[
U_i\setminus \bigl(U_i\cap \operatorname{Int}(M)\bigr)
=
U_i\cap \partial M.
\] Under the chart $\phi_i$, this set corresponds to the flat face of a half-ball, which has zero $n$-dimensional Lebesgue measure. Since Riemannian measure in coordinates is given by $\sqrt{\det(\mathbf{g})}\,d\lambda_n$, with a smooth strictly positive factor, it follows that $U_i\cap \partial M$ also has zero Riemannian measure. Therefore, \[
u=c_i
\quad\text{almost everywhere on }U_i.
\]

We now show that the local constants agree on overlaps. Let $i,j\in I$ be such that $U_i\cap U_j\neq\varnothing$. Since $U_i\cap U_j$ is a nonempty open subset of the manifold $M_\alpha$, there is $x\in U_i\cap U_j$. If $x\in \operatorname{Int}(M)$, some open neighborhood of $x$ in $M$ lies in $\operatorname{Int}(M)$, and hence \[
(U_i\cap U_j)\cap \operatorname{Int}(M)\neq\varnothing.
\] If $x\in\partial M$, since $U_i\cap U_j$ is open in the manifold with boundary, its expression in half-ball coordinates contains a small half-ball; in particular, it also contains points of $\operatorname{Int}(M)$. In either case, \[
(U_i\cap U_j)\cap \operatorname{Int}(M)
\] is a nonempty open subset of $\operatorname{Int}(M)$ and therefore has positive Riemannian measure.

But on that set we have, almost everywhere, \[
u=c_i
\qquad\text{and}\qquad
u=c_j.
\] It follows that $c_i=c_j$. Thus the constants $\{c_i\}_{i\in I}$ are compatible on nonempty intersections.

Since $M_\alpha$ is connected, we conclude that all constants $c_i$ agree with a single constant, denoted by $c_\alpha$. Moreover, since the cover is countable, the union of the exceptional null sets where the local equalities may fail still has measure zero. Therefore, \[
u=c_\alpha
\quad\text{almost everywhere on }M_\alpha.
\]

Since the connected component $M_\alpha$ was arbitrary, we conclude that $u$ is constant almost everywhere on each connected component of $M$. In particular, if $M$ is connected, there is $c\in\mathbb R$ such that \[
u=c
\quad\text{almost everywhere on }M.
\] \end{proof}

\begin{proposition}[The space of parallel sections and its projection] \label{prop:proyeccion-secciones-paralelas}\index{space of parallel sections and its projection@space of parallel sections and its projection} Let $(M,\mathbf{g})$ be a compact Riemannian manifold of dimension $n$ with or without boundary, and let $\mathbf{E}\to M$ be a smooth real vector bundle of finite rank equipped with a bundle metric $\mathbf{h}_{\mathbf{E}}$ and a compatible connection $\nabla^{\mathbf{E}}$. For fixed $1\leq p<\infty$, define $\displaystyle \mathcal K(\mathbf{E}):=\left\{\mathbf{v}\in W^{1,p}(M,\mathbf{E})\middle|\nabla_w^{\mathbf{E}}\mathbf{v}=0\right\}.$ Then this definition is independent of the exponent $p$, the space $\mathcal K(\mathbf{E})$ is finite-dimensional, and there is a continuous linear projection $\mathsf P_{\mathbf{E}}\colon L^p(M,\mathbf{E})\longrightarrow\mathcal K(\mathbf{E})$.

More precisely, if $(\boldsymbol{\sigma}_1,\dots,\boldsymbol{\sigma}_N)$ is an orthonormal basis of $\mathcal K(\mathbf{E})$ with respect to the inner product of $L^2(M,\mathbf{E})$, one may take $\displaystyle \mathsf P_{\mathbf{E}}\mathbf{u}:=\displaystyle\sum_{j=1}^N\left(\int_M\langle \mathbf{u},\boldsymbol{\sigma}_j\rangle_{\mathbf{h}_{\mathbf{E}}}\,d\lambda_{\mathbf{g}}\right)\boldsymbol{\sigma}_j,$ with the convention $\mathsf P_{\mathbf{E}}=0$ when $\mathcal K(\mathbf{E})=\{0\}$. \end{proposition}

\begin{proof} If $\mathbf{v}\in\mathcal K(\mathbf{E})$, the weak Kato inequality, Corollary~\ref{desigualdad debil de kato general}, implies that $|\mathbf{v}|_{\mathbf{h}_{\mathbf{E}}}\in W^{1,p}(M)$ and $\displaystyle \left|\nabla_w|\mathbf{v}|_{\mathbf{h}_{\mathbf{E}}}\right|_{\mathbf{g}}\leq|\nabla_w^{\mathbf{E}}\mathbf{v}|_{\mathbf{g},\mathbf{h}_{\mathbf{E}}}=0$ almost everywhere. Lemma~\ref{lem:derivada-debil-nula-funcion-constante} shows that $|\mathbf{v}|_{\mathbf{h}_{\mathbf{E}}}$ is constant almost everywhere on each connected component of $M$. Since $M$ is compact and has finitely many connected components, it follows that $\mathbf{v}\in L^\infty(M,\mathbf{E})$. Consequently, $\mathbf{v}\in L^q(M,\mathbf{E})$ for every $1\leq q<\infty$ and, since $\nabla_w^{\mathbf{E}}\mathbf{v}=0$ is independent of the integrability exponent, we have $\mathbf{v}\in W^{1,q}(M,\mathbf{E})$ and $\nabla_w^{\mathbf{E}}\mathbf{v}=0$ for every $q$. This proves that the kernel under consideration is independent of $p$.

Consider the closed ball $\displaystyle B:=\left\{\mathbf{v}\in\mathcal K(\mathbf{E})\middle|\|\mathbf{v}\|_{L^p(M,\mathbf{E})}\leq1\right\}$. Since $\nabla_w^{\mathbf{E}}\mathbf{v}=0$ for every $\mathbf{v}\in\mathcal K(\mathbf{E})$, the norms $L^p(M,\mathbf{E})$ and $W^{1,p}(M,\mathbf{E})$ agree on $\mathcal K(\mathbf{E})$. By Corollary~\ref{cor:rellich-descenso-un-orden-haces}, applied with $m=1$, $B$ is relatively compact in $L^p(M,\mathbf{E})$. Moreover, $B$ is closed in $L^p(M,\mathbf{E})$: if $\mathbf{v}_k\in B$ and $\mathbf{v}_k\to \mathbf{v}$ in $L^p(M,\mathbf{E})$, then $\nabla_w^{\mathbf{E}}\mathbf{v}_k=0\to0$ in $L^p(M,T^*M\otimes \mathbf{E})$ and, by Lemma~\ref{derivada debil es operador cerrado}, $\mathbf{v}\in W^{1,p}(M,\mathbf{E})$ and $\nabla_w^{\mathbf{E}}\mathbf{v}=0$. Thus $B$ is compact, and consequently $\mathcal K(\mathbf{E})$ is finite-dimensional.

Every element of $\mathcal K(\mathbf{E})$ belongs to $L^\infty(M,\mathbf{E})$ and, in particular, to $L^2(M,\mathbf{E})$, so a basis $(\boldsymbol{\sigma}_1,\dots,\boldsymbol{\sigma}_N)$ orthonormal with respect to the inner product of $L^2(M,\mathbf{E})$ may be chosen. For $\mathbf{u}\in L^p(M,\mathbf{E})$, the integrals in the definition of $\mathsf P_{\mathbf{E}}\mathbf{u}$ are well defined and, denoting the conjugate exponent of $p$ by $p'$, with $p'=\infty$ if $p=1$, Hölder's inequality (Proposition~\ref{desigualdad de holder}) gives $\displaystyle \|\mathsf P_{\mathbf{E}}\mathbf{u}\|_{L^p(M,\mathbf{E})}\leq\displaystyle\sum_{j=1}^N\|\boldsymbol{\sigma}_j\|_{L^p(M,\mathbf{E})}\|\boldsymbol{\sigma}_j\|_{L^{p'}(M,\mathbf{E})}\|\mathbf{u}\|_{L^p(M,\mathbf{E})}.$ Thus $\mathsf P_{\mathbf{E}}$ is linear and continuous. Finally, if $\mathbf{v}\in\mathcal K(\mathbf{E})$, its expansion in the orthonormal basis $(\boldsymbol{\sigma}_1,\dots,\boldsymbol{\sigma}_N)$ shows that $\mathsf P_{\mathbf{E}}\mathbf{v}=\mathbf{v}$, so $\mathsf P_{\mathbf{E}}$ is a projection onto $\mathcal K(\mathbf{E})$. \end{proof}

The Poincaré inequality can now be formulated as an estimate for a section modulo its parallel component.

\begin{theorem}[Poincaré inequality modulo the kernel of the connection] \label{teo:poincare-modulo-nucleo-conexion}\index{Poincare inequality@Poincaré inequality!modulo the kernel of the connection} Let $(M,\mathbf{g})$ be a compact Riemannian manifold of dimension $n$ with or without boundary, and let $\mathbf{E}\to M$ be a smooth real vector bundle of finite rank equipped with a bundle metric $\mathbf{h}_{\mathbf{E}}$ and a compatible connection $\nabla^{\mathbf{E}}$. If $1\leq p<\infty$, then there is a constant $C_P>0$, depending only on $M$, $\mathbf{g}$, $\mathbf{E}$, $\mathbf{h}_{\mathbf{E}}$, $\nabla^{\mathbf{E}}$, $p$, and the projection $\mathsf P_{\mathbf{E}}$, such that $\displaystyle \|\mathbf{u}-\mathsf P_{\mathbf{E}}\mathbf{u}\|_{L^p(M,\mathbf{E})}\leq C_P\|\nabla_w^{\mathbf{E}}\mathbf{u}\|_{L^p(M,T^*M\otimes \mathbf{E})}$ for every $\mathbf{u}\in W^{1,p}(M,\mathbf{E})$.

Consequently, $\displaystyle \|\mathbf{u}\|_{\nabla^{\mathbf{E}},p}:=\|\nabla_w^{\mathbf{E}}\mathbf{u}\|_{L^p(M,T^*M\otimes \mathbf{E})}$ defines a norm on $\ker(\mathsf P_{\mathbf{E}})\cap W^{1,p}(M,\mathbf{E})$ equivalent to the usual norm of $W^{1,p}(M,\mathbf{E})$. \end{theorem}

\begin{proof} Suppose, for a contradiction, that the inequality is false. Then, for each $k\in\mathbb N$, there is $\mathbf{u}_k\in W^{1,p}(M,\mathbf{E})$ such that $\displaystyle \|\mathbf{u}_k-\mathsf P_{\mathbf{E}}\mathbf{u}_k\|_{L^p(M,\mathbf{E})}>k\|\nabla_w^{\mathbf{E}}\mathbf{u}_k\|_{L^p(M,T^*M\otimes \mathbf{E})}$. Define $\displaystyle \mathbf{v}_k:=\displaystyle\frac{\mathbf{u}_k-\mathsf P_{\mathbf{E}}\mathbf{u}_k}{\|\mathbf{u}_k-\mathsf P_{\mathbf{E}}\mathbf{u}_k\|_{L^p(M,\mathbf{E})}}.$ Since $\mathsf P_{\mathbf{E}}^2=\mathsf P_{\mathbf{E}}$ and $\nabla_w^{\mathbf{E}}(\mathsf P_{\mathbf{E}}\mathbf{u}_k)=0$, we have $\mathsf P_{\mathbf{E}}\mathbf{v}_k=0$, $\|\mathbf{v}_k\|_{L^p(M,\mathbf{E})}=1$, and $\displaystyle \|\nabla_w^{\mathbf{E}}\mathbf{v}_k\|_{L^p(M,T^*M\otimes \mathbf{E})}<\displaystyle\frac{1}{k}$. In particular, $(\mathbf{v}_k)$ is bounded in $W^{1,p}(M,\mathbf{E})$.

By Corollary~\ref{cor:rellich-descenso-un-orden-haces}, applied with $m=1$, after passing to a subsequence there is $\mathbf{v}\in L^p(M,\mathbf{E})$ such that $\mathbf{v}_k\to \mathbf{v}$ in $L^p(M,\mathbf{E})$. On the other hand, $\nabla_w^{\mathbf{E}}\mathbf{v}_k\to0$ in $L^p(M,T^*M\otimes \mathbf{E})$ and, by Lemma~\ref{derivada debil es operador cerrado}, we have $\mathbf{v}\in W^{1,p}(M,\mathbf{E})$ and $\nabla_w^{\mathbf{E}}\mathbf{v}=0$. Therefore, $\mathbf{v}\in\mathcal K(\mathbf{E})$.

Continuity of $\mathsf P_{\mathbf{E}}$ on $L^p(M,\mathbf{E})$ and the identity $\mathsf P_{\mathbf{E}}\mathbf{v}_k=0$ imply $\mathsf P_{\mathbf{E}}\mathbf{v}=0$. However, since $\mathbf{v}\in\mathcal K(\mathbf{E})$ and $\mathsf P_{\mathbf{E}}$ is the identity on $\mathcal K(\mathbf{E})$, we also have $\mathsf P_{\mathbf{E}}\mathbf{v}=\mathbf{v}$. Consequently, $\mathbf{v}=0$, contradicting $\displaystyle \|\mathbf{v}\|_{L^p(M,\mathbf{E})}=\lim_{k\to\infty}\|\mathbf{v}_k\|_{L^p(M,\mathbf{E})}=1$. This proves the inequality.

If $\mathbf{u}\in\ker(\mathsf P_{\mathbf{E}})\cap W^{1,p}(M,\mathbf{E})$, then $\mathbf{u}-\mathsf P_{\mathbf{E}}\mathbf{u}=\mathbf{u}$ and therefore $\displaystyle \|\mathbf{u}\|_{L^p(M,\mathbf{E})}\leq C_P\|\nabla_w^{\mathbf{E}}\mathbf{u}\|_{L^p(M,T^*M\otimes \mathbf{E})}$. This gives $\displaystyle \|\nabla_w^{\mathbf{E}}\mathbf{u}\|_{L^p(M,T^*M\otimes \mathbf{E})}\leq\|\mathbf{u}\|_{W^{1,p}(M,\mathbf{E})}\leq(1+C_P^p)^{\frac{1}{p}}\|\nabla_w^{\mathbf{E}}\mathbf{u}\|_{L^p(M,T^*M\otimes \mathbf{E})},$ proving equivalence of norms. \end{proof}

\begin{corollary}[Poincaré inequality with the mean] \label{cor:poincare-promedio}\index{Poincare inequality with the mean@Poincaré inequality with the mean} Let $(M,\mathbf{g})$ be a compact connected Riemannian manifold of dimension $n$ with or without boundary, and let $1\leq p<\infty$. Then there is a constant $C_P>0$ such that $\displaystyle \|u-\overline u\|_{L^p(M)}\leq C_P\|\nabla_wu\|_{L^p(M,T^*M)}$ for every $u\in W^{1,p}(M)$, where $\displaystyle \overline u:=\displaystyle\frac{1}{\lambda_{\mathbf{g}}(M)}\int_Mu\,d\lambda_{\mathbf{g}}$. \end{corollary}

\begin{proof} For the trivial bundle $M\times\mathbb R$ equipped with the trivial connection, the space $\mathcal K(M\times\mathbb R)$ consists of constant functions and has dimension $1$. Since the constant function $1$ satisfies $\|1\|_{L^2(M)}=\lambda_{\mathbf{g}}(M)^{\frac{1}{2}}$, the function $\sigma=\lambda_{\mathbf{g}}(M)^{-\frac{1}{2}}$ has unit norm in $L^2(M)$, so $\{\sigma\}$ is an orthonormal basis of $\mathcal K(M\times\mathbb R)$. Consequently, $\displaystyle \mathsf P_{M\times\mathbb R}u=\left(\int_Mu\,\lambda_{\mathbf{g}}(M)^{-\frac{1}{2}}\,d\lambda_{\mathbf{g}}\right)\lambda_{\mathbf{g}}(M)^{-\frac{1}{2}}=\displaystyle\frac{1}{\lambda_{\mathbf{g}}(M)}\int_Mu\,d\lambda_{\mathbf{g}}=\overline u.$ The result follows from Theorem~\ref{teo:poincare-modulo-nucleo-conexion}. \end{proof} \begin{remark}\label{obs:poincare-compacto-no-conexa-sus-componentes-conexas-nucleo} If $M$ is disconnected and $M_1,\dots,M_\ell$ are its connected components, the kernel of the weak derivative consists of functions constant on each component. Consequently, the projection onto this kernel is given by $\mathsf P_{M\times\mathbb R}u
=
\displaystyle\sum_{i=1}^{\ell}
\left(
\displaystyle\frac{1}{\lambda_{\mathbf{g}}(M_i)}
\int_{M_i}u\,d\lambda_{\mathbf{g}}
\right)\mathbf{1}_{M_i}$.

This expression corresponds to subtracting the mean on each connected component. In particular, the global mean gives the correct formulation only when $M$ is connected. \end{remark} There are geometric situations in which a vector bundle admits no nontrivial parallel sections. In such cases the kernel of the weak covariant derivative is trivial, and the Poincaré inequality takes a particularly simple form. \begin{corollary}[The case without parallel sections] \label{cor:poincare-sin-secciones-paralelas}\index{case without parallel sections} Let $(M,\mathbf{g})$ be a compact Riemannian manifold of dimension $n$ with or without boundary, and let $\mathbf{E}\to M$ be a smooth real vector bundle of finite rank with a bundle metric and a compatible connection. Let $1\leq p<\infty$. If $\mathcal K(\mathbf{E})=\{0\}$, then there is a constant $C_P>0$ such that $\displaystyle \|\mathbf{u}\|_{L^p(M,\mathbf{E})}\leq C_P\|\nabla_w^{\mathbf{E}}\mathbf{u}\|_{L^p(M,T^*M\otimes \mathbf{E})}$ for every $\mathbf{u}\in W^{1,p}(M,\mathbf{E})$. In particular, $\|\nabla_w^{\mathbf{E}}\mathbf{u}\|_{L^p(M,T^*M\otimes \mathbf{E})}$ defines a norm on $W^{1,p}(M,\mathbf{E})$ equivalent to the usual Sobolev norm. \end{corollary}

\begin{proof} If $\mathcal K(\mathbf{E})=\{0\}$, then $\mathsf P_{\mathbf{E}}=0$ and the result follows directly from Theorem~\ref{teo:poincare-modulo-nucleo-conexion}. \end{proof}

The preceding formulation removes the parallel component of a section by a projection. When the manifold has boundary, another natural possibility is to impose a homogeneous Dirichlet condition. To prove the corresponding inequality, we need to relate the trace of a section to the trace of its bundle norm.

\begin{lemma}[Compatibility of the trace with the bundle norm] \label{lem:compatibilidad-traza-norma-fibrada}\index{compatibility of the trace with the bundle norm} Let $(M,\mathbf{g})$ be a compact Riemannian manifold with nonempty boundary, and let $\mathbf{E}\to M$ be a smooth real vector bundle of finite rank equipped with a bundle metric $\mathbf{h}_{\mathbf{E}}$ and a compatible connection $\nabla^{\mathbf{E}}$. If $1\leq p<\infty$ and $\mathbf{u}\in W^{1,p}(M,\mathbf{E})$, then $|\mathbf{u}|_{\mathbf{h}_{\mathbf{E}}}\in W^{1,p}(M)$ and $\operatorname{Tr}(|\mathbf{u}|_{\mathbf{h}_{\mathbf{E}}})=|\boldsymbol{\operatorname{Tr}}(\mathbf{u})|_{\mathbf{h}_{\mathbf{E}}}$ almost everywhere on $\partial M$. \end{lemma}

\begin{proof} By Corollary~\ref{desigualdad debil de kato general}, we have $|\mathbf{u}|_{\mathbf{h}_{\mathbf{E}}}\in W^{1,p}(M)$. By the Meyers--Serrin theorem for vector bundles on manifolds with boundary, Theorem~\ref{teo:meyers-serrin-haz-frontera}, there is a sequence $(\mathbf{u}_k)\subseteq\Gamma(\mathbf{E})\cap W^{1,p}(M,\mathbf{E})$ such that $\mathbf{u}_k\to \mathbf{u}$ in $W^{1,p}(M,\mathbf{E})$. The reverse triangle inequality implies $\bigl\||\mathbf{u}_k|_{\mathbf{h}_{\mathbf{E}}}-|\mathbf{u}|_{\mathbf{h}_{\mathbf{E}}}\bigr\|_{L^p(M)}\leq\|\mathbf{u}_k-\mathbf{u}\|_{L^p(M,\mathbf{E})}\to0$.

Applying Corollary~\ref{desigualdad debil de kato general} once more, now to each $\mathbf{u}_k$, gives $|\mathbf{u}_k|_{\mathbf{h}_{\mathbf{E}}}\in W^{1,p}(M)$ and \[
\left|\nabla_w|\mathbf{u}_k|_{\mathbf{h}_{\mathbf{E}}}\right|_{\mathbf{g}}
\leq
|\nabla_w^{\mathbf{E}} \mathbf{u}_k|_{\mathbf{g},\mathbf{h}_{\mathbf{E}}}
\] almost everywhere. Since $\nabla_w^{\mathbf{E}} \mathbf{u}_k\to\nabla_w^{\mathbf{E}} \mathbf{u}$ in $L^p(M,T^*M\otimes \mathbf{E})$, the sequence $\bigl(\nabla_w|\mathbf{u}_k|_{\mathbf{h}_{\mathbf{E}}}\bigr)$ is bounded in $L^p(M,T^*M)$.

If $1<p<\infty$, reflexivity of $L^p(M,T^*M)$ gives a subsequence and an element $\mathbf{Z}\in L^p(M,T^*M)$ such that \[
\nabla_w|\mathbf{u}_k|_{\mathbf{h}_{\mathbf{E}}}\rightharpoonup \mathbf{Z}
\qquad\text{in }L^p(M,T^*M).
\] If $p=1$, convergence $|\nabla_w^{\mathbf{E}} \mathbf{u}_k|_{\mathbf{g},\mathbf{h}_{\mathbf{E}}}\to|\nabla_w^{\mathbf{E}} \mathbf{u}|_{\mathbf{g},\mathbf{h}_{\mathbf{E}}}$ in $L^1(M)$ implies that the family \[
\left\{
|\nabla_w^{\mathbf{E}} \mathbf{u}_k|_{\mathbf{g},\mathbf{h}_{\mathbf{E}}}
\;\middle|\;
k\in\mathbb N
\right\}
\] is uniformly integrable. The weak Kato inequality then shows that \[
\left\{
\left|\nabla_w|\mathbf{u}_k|_{\mathbf{h}_{\mathbf{E}}}\right|_{\mathbf{g}}
\;\middle|\;
k\in\mathbb N
\right\}
\] is also uniformly integrable. By Theorem~\ref{dunford pettis para haces vectoriales}, applied to the bundle $T^*M$, the weak closure of $\{\nabla_w|\mathbf{u}_k|_{\mathbf{h}_{\mathbf{E}}}\mid k\in\mathbb N\}$ in $L^1(M,T^*M)$ is compact. The Eberlein--Šmulian theorem (Theorem~\ref{teo: eberlein smulian}) therefore allows us to extract a subsequence and find $\mathbf{Z}\in L^1(M,T^*M)$ such that \[
\nabla_w|\mathbf{u}_k|_{\mathbf{h}_{\mathbf{E}}}\rightharpoonup \mathbf{Z}
\qquad\text{in }L^1(M,T^*M).
\]

Thus in both cases, after passing to a subsequence, we may assume that $\nabla_w|\mathbf{u}_k|_{\mathbf{h}_{\mathbf{E}}}\rightharpoonup \mathbf{Z}$ in $L^p(M,T^*M)$ for some $\mathbf{Z}\in L^p(M,T^*M)$. For every $\boldsymbol{\eta}\in\Gamma_{c,\operatorname{Int}(M)}(T^*M)$, since $\nabla_w|\mathbf{u}_k|_{\mathbf{h}_{\mathbf{E}}}$ is the weak covariant derivative of $|\mathbf{u}_k|_{\mathbf{h}_{\mathbf{E}}}$, we have \[
\int_M
\left\langle\nabla_w|\mathbf{u}_k|_{\mathbf{h}_{\mathbf{E}}},\boldsymbol{\eta}\right\rangle_{\mathbf{g}}\,d\lambda_{\mathbf{g}}
=
\int_M|\mathbf{u}_k|_{\mathbf{h}_{\mathbf{E}}}\nabla^*\boldsymbol{\eta}\,d\lambda_{\mathbf{g}}.
\] Since $\boldsymbol{\eta}\in L^{p'}(M,T^*M)$ and $\nabla^*\boldsymbol{\eta}\in L^{p'}(M)$, with the convention $p'=\infty$ when $p=1$, weak convergence of $\nabla_w|\mathbf{u}_k|_{\mathbf{h}_{\mathbf{E}}}$ allows passage to the limit on the left-hand side. On the right-hand side we use strong convergence $|\mathbf{u}_k|_{\mathbf{h}_{\mathbf{E}}}\to|\mathbf{u}|_{\mathbf{h}_{\mathbf{E}}}$ in $L^p(M)$; indeed, by Hölder's inequality, \[
\left|
\int_M
\bigl(|\mathbf{u}_k|_{\mathbf{h}_{\mathbf{E}}}-|\mathbf{u}|_{\mathbf{h}_{\mathbf{E}}}\bigr)\nabla^*\boldsymbol{\eta}\,d\lambda_{\mathbf{g}}
\right|
\leq
\bigl\||\mathbf{u}_k|_{\mathbf{h}_{\mathbf{E}}}-|\mathbf{u}|_{\mathbf{h}_{\mathbf{E}}}\bigr\|_{L^p(M)}
\|\nabla^*\boldsymbol{\eta}\|_{L^{p'}(M)}
\longrightarrow0.
\] Consequently, \[
\int_M\langle \mathbf{Z},\boldsymbol{\eta}\rangle_{\mathbf{g}}\,d\lambda_{\mathbf{g}}
=
\int_M|\mathbf{u}|_{\mathbf{h}_{\mathbf{E}}}\nabla^*\boldsymbol{\eta}\,d\lambda_{\mathbf{g}},
\] so that $\mathbf{Z}=\nabla_w|\mathbf{u}|_{\mathbf{h}_{\mathbf{E}}}$.

Consider the isometric embedding \[
J_1:W^{1,p}(M)
\longrightarrow
L^p(M)\times L^p(M,T^*M),
\qquad
J_1(w):=(w,\nabla_w w),
\] where the product is equipped with the norm $\ell^p(\{1,2\})$, whose topology agrees with the product topology. Strong convergence $|\mathbf{u}_k|_{\mathbf{h}_{\mathbf{E}}}\to|\mathbf{u}|_{\mathbf{h}_{\mathbf{E}}}$ in $L^p(M)$ also implies weak convergence in that space and, together with $\nabla_w|\mathbf{u}_k|_{\mathbf{h}_{\mathbf{E}}}\rightharpoonup\nabla_w|\mathbf{u}|_{\mathbf{h}_{\mathbf{E}}}$ in $L^p(M,T^*M)$, Proposition~\ref{prop:topologia-debil-producto-finito} gives \[
J_1(|\mathbf{u}_k|_{\mathbf{h}_{\mathbf{E}}})
\rightharpoonup
J_1(|\mathbf{u}|_{\mathbf{h}_{\mathbf{E}}})
\] in $L^p(M)\times L^p(M,T^*M)$. By Theorem~\ref{encaje isometrico en Lp sobolev haces}, applied to the trivial bundle, the image of $J_1$ is closed. Proposition~\ref{prop:topologia-debil-subespacio-normado} then shows that this convergence is also weak in the image of $J_1$. Since $J_1\colon W^{1,p}(M)\longrightarrow\operatorname{Im}(J_1)$ is an isometric isomorphism, its inverse $J_1^{-1}\colon \operatorname{Im}(J_1)\longrightarrow W^{1,p}(M)$ is linear and continuous and, by Proposition~\ref{prop:operador-lineal-continuo-topologias-debiles}, is continuous for the weak topologies. Applying $J_1^{-1}$ to the preceding convergence, we conclude that \[
|\mathbf{u}_k|_{\mathbf{h}_{\mathbf{E}}}\rightharpoonup|\mathbf{u}|_{\mathbf{h}_{\mathbf{E}}}
\qquad\text{in }W^{1,p}(M).
\]

By Theorem~\ref{teo:traza-entera-haces}, applied respectively to the trivial bundle and to $\mathbf{E}$, the trace operators $\operatorname{Tr}\colon W^{1,p}(M)\longrightarrow L^p(\partial M)$ and $\boldsymbol{\operatorname{Tr}}\colon W^{1,p}(M,\mathbf{E})\longrightarrow L^p(\partial M,\mathbf{E}\restriction_{\partial M})$ are linear and continuous. By Proposition~\ref{prop:operador-lineal-continuo-topologias-debiles}, the operator $\operatorname{Tr}$ is continuous for the weak topologies, so $\operatorname{Tr}(|\mathbf{u}_k|_{\mathbf{h}_{\mathbf{E}}})\rightharpoonup\operatorname{Tr}(|\mathbf{u}|_{\mathbf{h}_{\mathbf{E}}})$ in $L^p(\partial M)$. On the other hand, continuity of $\boldsymbol{\operatorname{Tr}}$ and convergence $\mathbf{u}_k\to \mathbf{u}$ in $W^{1,p}(M,\mathbf{E})$ imply $\boldsymbol{\operatorname{Tr}}(\mathbf{u}_k)\to\boldsymbol{\operatorname{Tr}}(\mathbf{u})$ in $L^p(\partial M,\mathbf{E}\restriction_{\partial M})$. The reverse triangle inequality then gives \[
|\boldsymbol{\operatorname{Tr}}(\mathbf{u}_k)|_{\mathbf{h}_{\mathbf{E}}}
\longrightarrow
|\boldsymbol{\operatorname{Tr}}(\mathbf{u})|_{\mathbf{h}_{\mathbf{E}}}
\qquad\text{in }L^p(\partial M).
\]

To justify this identity for each smooth section $\mathbf{u}_k$, fix $k$ and use the smooth functions \[
f_{k,\varepsilon}:=\sqrt{|\mathbf{u}_k|_{\mathbf{h}_{\mathbf{E}}}^{\,2}+\varepsilon^2}-\varepsilon,
\qquad \varepsilon>0.
\] The proof of Lemma~\ref{lem:derivada-debil-modulo-seccion-suave}, applied on the compact set $M$, proves that $f_{k,\varepsilon}\to|\mathbf{u}_k|_{\mathbf{h}_{\mathbf{E}}}$ in $W^{1,p}(M)$ for every $1\leq p<\infty$: the functions converge in $L^p(M)$ and their differentials in $L^p(M,T^*M)$ by dominated convergence. Continuity of the trace then gives $\operatorname{Tr}(f_{k,\varepsilon})\to\operatorname{Tr}(|\mathbf{u}_k|_{\mathbf{h}_{\mathbf{E}}})$ in $L^p(\partial M)$. Since $f_{k,\varepsilon}$ is smooth, its trace is its restriction, and the inequality \[
0\leq|\mathbf{u}_k|_{\mathbf{h}_{\mathbf{E}}}-f_{k,\varepsilon}\leq\varepsilon
\] proves that these restrictions converge uniformly on $\partial M$ to $|\mathbf{u}_k\restriction_{\partial M}|_{\mathbf{h}_{\mathbf{E}}}$. Since $\partial M$ has finite measure, convergence also holds in $L^p(\partial M)$. Uniqueness of the limit gives $\operatorname{Tr}(|\mathbf{u}_k|_{\mathbf{h}_{\mathbf{E}}})=|\boldsymbol{\operatorname{Tr}}(\mathbf{u}_k)|_{\mathbf{h}_{\mathbf{E}}}$ almost everywhere on $\partial M$ for every $k$. The preceding strong convergence also implies weak convergence, so the same sequence $\operatorname{Tr}(|\mathbf{u}_k|_{\mathbf{h}_{\mathbf{E}}})=|\boldsymbol{\operatorname{Tr}}(\mathbf{u}_k)|_{\mathbf{h}_{\mathbf{E}}}$ converges weakly both to $\operatorname{Tr}(|\mathbf{u}|_{\mathbf{h}_{\mathbf{E}}})$ and to $|\boldsymbol{\operatorname{Tr}}(\mathbf{u})|_{\mathbf{h}_{\mathbf{E}}}$. Uniqueness of the weak limit finally gives \[
\operatorname{Tr}(|\mathbf{u}|_{\mathbf{h}_{\mathbf{E}}})
=
|\boldsymbol{\operatorname{Tr}}(\mathbf{u})|_{\mathbf{h}_{\mathbf{E}}}
\] almost everywhere on $\partial M$. \end{proof}

We can now prove the Poincaré inequality for sections satisfying a homogeneous Dirichlet condition.

\begin{theorem}[Poincaré inequality in $W^{1,p}_0(M,\mathbf{E})$] \label{teo:poincare-W0-haces}\index{Poincare inequality@Poincaré inequality} Let $(M,\mathbf{g})$ be a compact Riemannian manifold of dimension $n$ with nonempty smooth boundary, and suppose that every connected component of $M$ has nonempty boundary. Let $\mathbf{E}\to M$ be a smooth real vector bundle of finite rank equipped with a bundle metric $\mathbf{h}_{\mathbf{E}}$ and a compatible connection $\nabla^{\mathbf{E}}$. If $1\leq p<\infty$, then there is a constant $C_P>0$, depending only on $M$, $\mathbf{g}$, $\mathbf{E}$, $\mathbf{h}_{\mathbf{E}}$, $\nabla^{\mathbf{E}}$, and $p$, such that $\displaystyle \|\mathbf{u}\|_{L^p(M,\mathbf{E})}\leq C_P\|\nabla_w^{\mathbf{E}}\mathbf{u}\|_{L^p(M,T^*M\otimes \mathbf{E})}$ for every $\mathbf{u}\in W^{1,p}_0(M,\mathbf{E})$. Consequently, $\displaystyle \|\mathbf{u}\|_{\nabla^{\mathbf{E}},p}:=\left(\int_M|\nabla_w^{\mathbf{E}}\mathbf{u}|_{\mathbf{g},\mathbf{h}_{\mathbf{E}}}^p\,d\lambda_{\mathbf{g}}\right)^{\frac{1}{p}}$ defines a norm on $W^{1,p}_0(M,\mathbf{E})$ equivalent to the usual norm of $W^{1,p}(M,\mathbf{E})$. \end{theorem}

\begin{proof} Suppose, for a contradiction, that the inequality is false. Then there is a sequence $(\mathbf{u}_k)\subseteq W^{1,p}_0(M,\mathbf{E})$ such that $\|\mathbf{u}_k\|_{L^p(M,\mathbf{E})}=1$ for every $k\in\mathbb N$ and $\displaystyle \|\nabla_w^{\mathbf{E}}\mathbf{u}_k\|_{L^p(M,T^*M\otimes \mathbf{E})}\to0$. In particular, $(\mathbf{u}_k)$ is bounded in $W^{1,p}(M,\mathbf{E})$. By Corollary~\ref{cor:rellich-descenso-un-orden-haces}, applied with $m=1$, after passing to a subsequence there is $\mathbf{u}\in L^p(M,\mathbf{E})$ such that $\mathbf{u}_k\to \mathbf{u}$ in $L^p(M,\mathbf{E})$.

Since $\nabla_w^{\mathbf{E}}\mathbf{u}_k\to0$ in $L^p(M,T^*M\otimes \mathbf{E})$ and the weak covariant derivative is a closed operator from $L^p(M,\mathbf{E})$ to $L^p(M,T^*M\otimes \mathbf{E})$, we have $\mathbf{u}\in W^{1,p}(M,\mathbf{E})$ and $\nabla_w^{\mathbf{E}}\mathbf{u}=0$. Therefore, $\mathbf{u}_k\to \mathbf{u}$ in $W^{1,p}(M,\mathbf{E})$. This step uses only strong convergence and closedness of the weak derivative, Lemma~\ref{derivada debil es operador cerrado}, and is also valid for $p=1$; it does not require reflexivity. Since $W^{1,p}_0(M,\mathbf{E})$ is a closed subspace of $W^{1,p}(M,\mathbf{E})$, it follows that $\mathbf{u}\in W^{1,p}_0(M,\mathbf{E})$ and, by Theorem~\ref{teo:caracterizacion-W0-trazas} applied with $m=1$, $\boldsymbol{\operatorname{Tr}}(\mathbf{u})=0$.

Write $\displaystyle M=\bigsqcup_{i=1}^{\ell}M_i$, where $M_1,\dots,M_\ell$ are the connected components of $M$. By the weak Kato inequality, Corollary~\ref{desigualdad debil de kato general}, we have $|\mathbf{u}|_{\mathbf{h}_{\mathbf{E}}}\in W^{1,p}(M)$ and $\displaystyle \left|\nabla_w|\mathbf{u}|_{\mathbf{h}_{\mathbf{E}}}\right|_{\mathbf{g}}\leq|\nabla_w^{\mathbf{E}}\mathbf{u}|_{\mathbf{g},\mathbf{h}_{\mathbf{E}}}=0$ almost everywhere. Lemma~\ref{lem:derivada-debil-nula-funcion-constante} implies that, for each $i\in\{1,\dots,\ell\}$, there is a constant $c_i\geq0$ such that $|\mathbf{u}|_{\mathbf{h}_{\mathbf{E}}}=c_i$ almost everywhere on $M_i$.

By Lemma~\ref{lem:compatibilidad-traza-norma-fibrada}, $\operatorname{Tr}(|\mathbf{u}|_{\mathbf{h}_{\mathbf{E}}})=|\boldsymbol{\operatorname{Tr}}(\mathbf{u})|_{\mathbf{h}_{\mathbf{E}}}=0$ almost everywhere on $\partial M$. Since $\partial M_i\neq\varnothing$ and the trace of the constant function $c_i$ on $\partial M_i$ is $c_i$, we conclude that $c_i=0$ for each $i\in\{1,\dots,\ell\}$. Thus $|\mathbf{u}|_{\mathbf{h}_{\mathbf{E}}}=0$ almost everywhere on $M$ and consequently $\mathbf{u}=0$ almost everywhere.

However, convergence $\mathbf{u}_k\to \mathbf{u}$ in $L^p(M,\mathbf{E})$ implies $\displaystyle \|\mathbf{u}\|_{L^p(M,\mathbf{E})}=\lim_{k\to\infty}\|\mathbf{u}_k\|_{L^p(M,\mathbf{E})}=1$, a contradiction. Thus there is $C_P>0$ such that $\displaystyle \|\mathbf{u}\|_{L^p(M,\mathbf{E})}\leq C_P\|\nabla_w^{\mathbf{E}}\mathbf{u}\|_{L^p(M,T^*M\otimes \mathbf{E})}$ for every $\mathbf{u}\in W^{1,p}_0(M,\mathbf{E})$.

Finally, for every $\mathbf{u}\in W^{1,p}_0(M,\mathbf{E})$ we have $\displaystyle \|\nabla_w^{\mathbf{E}}\mathbf{u}\|_{L^p(M,T^*M\otimes \mathbf{E})}\leq\|\mathbf{u}\|_{W^{1,p}(M,\mathbf{E})}\leq(1+C_P^p)^{\frac{1}{p}}\|\nabla_w^{\mathbf{E}}\mathbf{u}\|_{L^p(M,T^*M\otimes \mathbf{E})}$, so $\|\cdot\|_{\nabla^{\mathbf{E}},p}$ is a norm equivalent to the usual norm of $W^{1,p}(M,\mathbf{E})$. \end{proof}
\chapter{Distributions on vector bundles} \label{cap:distribuciones-en-haces-vectoriales}

Distribution theory on a manifold must respect the geometry of the space while extending the duality used in the Euclidean case. For scalar functions, a distribution can be viewed as a continuous functional on test functions. On a vector bundle, the natural formulation requires pairing test sections of the dual bundle and controlling their support and covariant derivatives.

We first develop the topology of spaces of compactly supported smooth sections. We then use it to define bundle-valued distributions, their restrictions, their support, and the distributional action of differential operators. The formal adjoint allows an operator to be transferred to the test section, just as integration by parts defines distributional derivatives on $\mathbb R^n$.

This framework will be needed later to introduce Sobolev spaces of real order and operators defined by functional calculus, since it allows their elements and equations to be interpreted without recourse to particular coordinates.

\begin{semblanzaHistorica}{Duality as a defining principle} Schwartz's theory taught us that a singular object can be understood through its action on test functions. On a vector bundle, this idea requires care concerning which fibers are paired, which measure permits integration, and how an operator is transferred to the test factor. The result is a theory that preserves the local nature of distributions while respecting the global geometry of the bundle. This duality will underlie elliptic regularity and the kernel theorem appearing at the end of the chapter~\cite{Schwartz1950}. \end{semblanzaHistorica}

\section{The space of test sections}\label{sec:espacio-secciones-prueba}

In this section, we construct the spaces of test sections associated with a smooth vector bundle. The idea parallels the construction of $\mathcal D(\Omega)$ in the Euclidean case, but the $C^m$ norms are now defined using covariant derivatives and bundle metrics. This formulation allows us to work intrinsically while retaining the functional properties needed to define bundle-valued distributions.

For a boundary chart, its image $\Omega$ is a relatively open subset of the closed half-space. The notation $C^\infty(\Omega)$ means smoothness up to the boundary, and $\mathcal D(\Omega)$ allows compact supports meeting it. Its dual is understood with this convention. All interiors and neighborhoods of compact subsets of a manifold are taken in its relative topology. The extension by zero of a compactly supported function from a relatively open subset to the half-space is smooth; to extend it across the geometric boundary, we use Lemma~\ref{lem:extension-suave-semiespacio-soporte-fijo}.

\begin{definition}\label{def:normas-Cm-secciones-soporte-compacto}\index{section!\(C^m\) norm on a compact set} Let $(M,\mathbf{g})$ be a Riemannian manifold with or without boundary, let $\mathbf{E}\longrightarrow M$ be a smooth vector bundle equipped with a bundle metric $\mathbf{h}_{\mathbf{E}}$, let $\nabla^{TM}$ be a connection on $TM$, and let $\nabla^{\mathbf{E}}$ be a connection on $\mathbf{E}$. Let $K\subseteq M$ be compact. For each $\boldsymbol{\phi}\in \Gamma^{m}(\mathbf{E})$, define \[
\|\boldsymbol{\phi}\|_{C^{m}(K)}:=\max_{0\le j\le m}\sup_{x\in K}|\nabla^{j}\boldsymbol{\phi}(x)|_{\mathbf{g},\mathbf{h}_{\mathbf{E}}},
\] where $\nabla^{j}\boldsymbol{\phi}$ denotes the covariant derivative of order $j$ of $\boldsymbol{\phi}$, regarded as a section of the tensor bundle $T^{(0,j)}(TM)\otimes \mathbf{E}$ and induced by the connections $\nabla^{TM}$ and $\nabla^{\mathbf{E}}$ as in Definition~\ref{derivada covariante de orden superior para haces vectoriales}. The norm $|\cdot|_{\mathbf{g},\mathbf{h}_{\mathbf{E}}}$ is the norm induced by $\mathbf{g}$ and $\mathbf{h}_{\mathbf{E}}$ on the corresponding tensor bundles; see Proposition~\ref{metrica en tensores con E general}.

Also define \[
\Gamma_{K}^{m}(\mathbf{E}):=\{\boldsymbol{\phi}\in \Gamma^{m}(\mathbf{E})\mid \operatorname{supp}(\boldsymbol{\phi})\subseteq K\}.
\] When $m=\infty$, we write \[
\Gamma_K(\mathbf{E}):=\{\boldsymbol{\phi}\in \Gamma(\mathbf{E})\mid \operatorname{supp}(\boldsymbol{\phi})\subseteq K\}.
\] \end{definition}

For each $m\in\mathbb N_0$, the quantity $\|\cdot\|_{C^m(K)}$ defines a norm on $\Gamma_K^m(\mathbf{E})$. The triangle inequality and homogeneity are immediate. If $\|\phi\|_{C^m(K)}=0$, then in particular $\phi=0$ on $K$. Since $\operatorname{supp}(\phi)\subseteq K$, it follows that $\phi=0$ on $M\setminus K$, and hence $\phi=0$ throughout $M$.

Although the preceding definition uses a Riemannian metric, a bundle metric, and auxiliary connections, the resulting topology on $\Gamma_K^m(\mathbf{E})$ is independent of these choices. This matters because the subsequent constructions should depend on the bundle and manifold, rather than on particular auxiliary data.

\begin{theorem}\label{norma ck independiente de metrica y conexion}\index{norm on sections!independence of auxiliary metrics and connections} Let $M$ be a smooth manifold with or without boundary, let $\mathbf{g}$ and $\widetilde{\mathbf{g}}$ be two Riemannian metrics on $M$, let $\nabla^{TM}$ and $\widetilde{\nabla}^{TM}$ be two connections on $TM$, let $\mathbf{E}\longrightarrow M$ be a smooth vector bundle of finite rank, let $\mathbf{h}_{\mathbf{E}}$ and $\widetilde{\mathbf{h}}_{\mathbf{E}}$ be two smooth bundle metrics on $\mathbf{E}$, and let $\nabla^{\mathbf{E}}$ and $\widetilde{\nabla}^{\mathbf{E}}$ be two connections on $\mathbf{E}$. Also, let $K\subseteq M$ be compact. Then the norms $\|\cdot\|_{C^m(K)}$ and $\|\cdot\|_{\widetilde{C^m}(K)}$ are equivalent on $\Gamma_K^m(\mathbf{E})$. \end{theorem}

\begin{proof} For $\mathbf{u}\in\Gamma_K^m(\mathbf{E})$, write \[
\|\mathbf{u}\|_{C^m(K)}
=
\max_{0\leq j\leq m}\sup_{x\in K}|\nabla^j\mathbf{u}(x)|_{\mathbf{g},\mathbf{h}_{\mathbf{E}}},
\qquad
\|\mathbf{u}\|_{\widetilde{C^m}(K)}
=
\max_{0\leq j\leq m}\sup_{x\in K}|\widetilde{\nabla}^j\mathbf{u}(x)|_{\widetilde{\mathbf{g}},\widetilde{\mathbf{h}}_{\mathbf{E}}}.
\] Applying Lemma~\ref{lema: independencia de normas en compactos} with the integer $m$ gives a constant $C_m>0$ such that, for every $x\in K$, \[
\displaystyle\sum_{j=0}^{m}|\widetilde{\nabla}^{j}\mathbf{u}(x)|^{2}_{\widetilde{\mathbf{g}},\widetilde{\mathbf{h}}_{\mathbf{E}}}
\leq
C_m\displaystyle\sum_{j=0}^{m}|\nabla^{j}\mathbf{u}(x)|^{2}_{\mathbf{g},\mathbf{h}_{\mathbf{E}}}.
\] Therefore \[
\max_{0\leq j\leq m}|\widetilde{\nabla}^{j}\mathbf{u}(x)|^{2}_{\widetilde{\mathbf{g}},\widetilde{\mathbf{h}}_{\mathbf{E}}}
\leq
C_m(m+1)\max_{0\leq j\leq m}|\nabla^{j}\mathbf{u}(x)|^{2}_{\mathbf{g},\mathbf{h}_{\mathbf{E}}}.
\] Taking the supremum over $x\in K$, we obtain \[
\|\mathbf{u}\|_{\widetilde{C^m}(K)}
\leq
\sqrt{C_m(m+1)}\,\|\mathbf{u}\|_{C^m(K)}.
\] The reverse argument, interchanging the data with and without tildes, produces a constant $\widetilde C>0$ such that \[
\|\mathbf{u}\|_{C^m(K)}
\leq
\widetilde C\,\|\mathbf{u}\|_{\widetilde{C^m}(K)}.
\] These two inequalities give the asserted comparison constants. \end{proof}

To prove that $\Gamma_K^m(\mathbf{E})$ is complete, we adapt the Euclidean argument to uniform limits of sections and their covariant derivatives. This intrinsic formulation will later allow passage to the limit in covariant identities.

\begin{lemma}[Uniform limit of continuous sections]\label{lema: limite uniforme de secciones continuas}\index{uniform limit of continuous sections@uniform limit of continuous sections} Let $M$ be a smooth manifold with or without boundary, and let $\mathbf{F}\longrightarrow M$ be a vector bundle of finite rank equipped with a continuous bundle metric $\mathbf{h}_{\mathbf{F}}$. Let $(\mathbf{s}_n)_{n\in\mathbb N}$ be a sequence of continuous sections of $\mathbf{F}$, and suppose that there exists a section $\mathbf{s}\colon M\longrightarrow \mathbf{F}$ such that \[
\sup_{x\in M}|\mathbf s_n(x)-\mathbf s(x)|_{\mathbf h_{\mathbf F}(x)}\longrightarrow 0.
\] Then $\mathbf{s}$ is continuous. \end{lemma}

\begin{proof} Fix $p\in M$. Choose a local trivialization $\Psi\colon \mathbf{F}\restriction_U\longrightarrow U\times\mathbb K^r$ over an open neighborhood $U$ of $p$, where $r$ is the rank of $\mathbf{F}$ and $\mathbb K=\mathbb R$ or $\mathbb C$. Also choose an open neighborhood $V$ of $p$ such that $\overline V\subseteq U$ and $\overline V$ is compact.

For each $n\in\mathbb N$, write $\Psi(\mathbf{s}_n(x))=(x,u_n(x))$ and $\Psi(\mathbf{s}(x))=(x,u(x))$ for $x\in U$, where $u_n,u\colon U\longrightarrow\mathbb K^r$.

Since the bundle metric $\mathbf{h}_{\mathbf{F}}$ is continuous and $\overline V$ is compact, the norms induced by $\mathbf{h}_{\mathbf{F}}$ on the fibers of $\mathbf{F}\restriction_{\overline V}$ are uniformly equivalent to the Euclidean norm on $\mathbb K^r$. Thus there exists a constant $C>0$ such that \[
|u_n(x)-u(x)|
\leq
C\,|\mathbf{s}_n(x)-\mathbf{s}(x)|_{\mathbf{h}_{\mathbf{F}}}
\] for every $x\in\overline V$ and every $n\in\mathbb N$.

Consequently, \[
\|u_n-u\|_{C^0(\overline V,\mathbb K^r)}
\leq
C\sup_{x\in M}|\mathbf s_n(x)-\mathbf s(x)|_{\mathbf h_{\mathbf F}(x)}
\longrightarrow 0.
\] Thus $u_n\to u$ uniformly on $\overline V$. Since each $u_n$ is continuous, $u$ is also continuous on $V$. Therefore $\mathbf{s}$ is continuous on $V$.

Since $p\in M$ was arbitrary, we conclude that $\mathbf{s}$ is continuous throughout $M$. \end{proof}

\begin{proposition}\label{prop: GammaKm es Banach} Let $(M,\mathbf{g})$ be a Riemannian manifold with or without boundary, let $\mathbf{E}\longrightarrow M$ be a smooth vector bundle equipped with a bundle metric $\mathbf{h}_{\mathbf{E}}$, let $\nabla^{TM}$ be a connection on $TM$, let $\nabla^{\mathbf{E}}$ be a connection on $\mathbf{E}$, and let $K\subseteq M$ be compact. For each $m\in\mathbb N_0$, the space $\bigl(\Gamma_K^m(\mathbf{E}),\|\cdot\|_{C^m(K)}\bigr)$ is a Banach space. \end{proposition}

\begin{proof}
For a continuous section $\mathbf a$, the essential supremum of $|\mathbf a|$ equals its pointwise supremum: Riemannian measure assigns positive measure to every nonempty relatively open set, including in boundary charts, and $|\mathbf a|$ is continuous. We use this equality for differences of continuous sections. For the limit, whose continuity is not yet known, we prove pointwise uniform convergence directly.

 We use the convention $\nabla^0\boldsymbol{\phi}=\boldsymbol{\phi}$. For each $r\in\{0,\dots,m\}$, denote by $\mathbf{F}_r:=T^{(0,r)}(TM)\otimes \mathbf{E}$ the tensor bundle in which $\nabla^r\boldsymbol{\phi}$ takes values. This bundle is equipped with the connection induced by $\nabla^{TM}$ and $\nabla^{\mathbf{E}}$, and the bundle metric induced by $\mathbf{g}$ and $\mathbf{h}_{\mathbf{E}}$.

Let $(\boldsymbol{\phi}_n)_{n\in\mathbb N}$ be a Cauchy sequence in $\Gamma_K^m(\mathbf{E})$. Fix $r\in\{0,\dots,m\}$.

If $x\in M\setminus K$, since $\operatorname{supp}(\boldsymbol{\phi}_n)\subseteq K$, there is an open neighborhood of $x$ on which $\boldsymbol{\phi}_n$ vanishes. Thus $\nabla^r\boldsymbol{\phi}_n(x)=0$. Consequently, $\nabla^r\boldsymbol{\phi}_n$ vanishes on $M\setminus K$ for every $n\in\mathbb N$.

Moreover, for any $n,p\in\mathbb N$, \[
\|\nabla^r\boldsymbol{\phi}_n-\nabla^r\boldsymbol{\phi}_p\|_{L^\infty(M,T^{(0,r)}(TM)\otimes\mathbf E)}
=
\|\nabla^r\boldsymbol{\phi}_n-\nabla^r\boldsymbol{\phi}_p\|_{C^0(K,T^{(0,r)}(TM)\otimes\mathbf E)}
\leq
\|\boldsymbol{\phi}_n-\boldsymbol{\phi}_p\|_{C^m(K)}.
\] Since $(\boldsymbol{\phi}_n)_{n\in\mathbb N}$ is Cauchy in $\Gamma_K^m(\mathbf{E})$, it follows that $(\nabla^r\boldsymbol{\phi}_n)_{n\in\mathbb N}$ is uniformly Cauchy on $M$.

For each $x\in M$, the sequence $(\nabla^r\boldsymbol{\phi}_n(x))_{n\in\mathbb N}$ is Cauchy in the fiber $(\mathbf{F}_r)_x$, which is a finite-dimensional normed vector space and hence complete. Define \[
\mathbf{A}_r(x):=
\lim_{n\to\infty}\nabla^r\boldsymbol{\phi}_n(x),
\qquad x\in M.
\] This gives a section $\mathbf{A}_r\colon M\longrightarrow \mathbf{F}_r$.

Let us show that convergence is uniform. Let $\varepsilon>0$. Since $(\nabla^r\boldsymbol{\phi}_n)_{n\in\mathbb N}$ is uniformly Cauchy, there exists $N\in\mathbb N$ such that, if $n,p\geq N$, then \[
\|\nabla^r\boldsymbol{\phi}_n-\nabla^r\boldsymbol{\phi}_p\|_{L^\infty(M,T^{(0,r)}(TM)\otimes\mathbf E)}
<
\varepsilon.
\] Fix $n\geq N$ and let $p\to\infty$. For each $x\in M$, \[
\left|
\nabla^r\boldsymbol{\phi}_n(x)-\mathbf{A}_r(x)
\right|_{\mathbf{g},\mathbf{h}_{\mathbf{E}}}
\leq
\varepsilon.
\] Taking the supremum over $x\in M$ gives \[
\sup_{x\in M}|\nabla^r\boldsymbol{\phi}_n(x)-\mathbf A_r(x)|_{\mathbf g,\mathbf h_{\mathbf E}}
\leq
\varepsilon.
\] Therefore $\nabla^r\boldsymbol{\phi}_n\to \mathbf{A}_r$ uniformly on $M$.

Since each $\nabla^r\boldsymbol{\phi}_n$ is continuous, Lemma~\ref{lema: limite uniforme de secciones continuas} implies that $\mathbf{A}_r$ is continuous. Moreover, since $\nabla^r\boldsymbol{\phi}_n=0$ on $M\setminus K$ for every $n$, passage to the limit gives $\mathbf{A}_r=0$ on $M\setminus K$.

Define $\boldsymbol{\phi}:=\mathbf{A}_0$. Then $\boldsymbol{\phi}$ is a continuous section of $\mathbf{E}$ and vanishes on $M\setminus K$. Since $K$ is compact and hence closed, it follows that $\operatorname{supp}(\boldsymbol{\phi})\subseteq K$.

It remains to prove that $\boldsymbol{\phi}$ is of class $C^m$ and that $\nabla^r\boldsymbol{\phi}=\mathbf{A}_r$ for every $r\in\{0,\dots,m\}$.

Fix $r\in\{0,\dots,m-1\}$. Consider the bundle $\mathbf{F}_r=T^{(0,r)}(TM)\otimes \mathbf{E}$ with the induced connection. Let $\gamma\colon [a,b]\longrightarrow M$ be a smooth curve, and let $t,t_0\in[a,b]$.

Since $\boldsymbol{\phi}_n\in\Gamma_K^m(\mathbf{E})$ and $r\leq m-1$, the section $\nabla^r\boldsymbol{\phi}_n$ is of class $C^1$ as a section of $\mathbf{F}_r$. Apply Lemma~\ref{lema: TFC covariante} to the section $\nabla^r\boldsymbol{\phi}_n$ restricted to $\gamma$. We obtain \[
P_{\gamma,t\to t_0}^{\mathbf{F}_r}
\bigl(
\nabla^r\boldsymbol{\phi}_n(\gamma(t))
\bigr)
-
\nabla^r\boldsymbol{\phi}_n(\gamma(t_0))
\] \[
\qquad =
\int_{t_0}^{t}
P_{\gamma,\tau\to t_0}^{\mathbf{F}_r}
\left(
\nabla^{\mathbf{F}_r}_{\dot\gamma(\tau)}
\bigl(\nabla^r\boldsymbol{\phi}_n\bigr)(\gamma(\tau))
\right)
\,d\tau.
\]

By the definition of the higher-order covariant derivative, under the identification $T^*M\otimes \mathbf{F}_r\cong T^{(0,r+1)}(TM)\otimes \mathbf{E}$ we have $\nabla^{\mathbf{F}_r}(\nabla^r\boldsymbol{\phi}_n)=\nabla^{r+1}\boldsymbol{\phi}_n$. Therefore \[
\nabla^{\mathbf{F}_r}_{\dot\gamma(\tau)}
\bigl(\nabla^r\boldsymbol{\phi}_n\bigr)(\gamma(\tau))
=
\iota^{\mathrm{prim}}_{\dot\gamma(\tau)}
\bigl(\nabla^{r+1}\boldsymbol{\phi}_n(\gamma(\tau))\bigr),
\] where $\iota^{\mathrm{prim}}_v$ denotes contraction in the first covariant factor; that is, \[
\bigl(\iota^{\mathrm{prim}}_vA\bigr)(X_1,\ldots,X_r)
=A(v,X_1,\ldots,X_r).
\] Thus \[
P_{\gamma,t\to t_0}^{\mathbf{F}_r}
\bigl(
\nabla^r\boldsymbol{\phi}_n(\gamma(t))
\bigr)
-
\nabla^r\boldsymbol{\phi}_n(\gamma(t_0))
\] \[
\qquad =
\int_{t_0}^{t}
P_{\gamma,\tau\to t_0}^{\mathbf{F}_r}
\left(
\iota^{\mathrm{prim}}_{\dot\gamma(\tau)}
\bigl(\nabla^{r+1}\boldsymbol{\phi}_n(\gamma(\tau))\bigr)
\right)
\,d\tau.
\]

We let $n\to\infty$ in this identity. Since $\nabla^r\boldsymbol{\phi}_n\to \mathbf{A}_r$ uniformly on $M$, the left-hand side converges to \[
P_{\gamma,t\to t_0}^{\mathbf{F}_r}
\bigl(
\mathbf{A}_r(\gamma(t))
\bigr)
-
\mathbf{A}_r(\gamma(t_0)).
\]

Let us now check that we can also pass to the limit in the integral. For fixed $t$ and $t_0$, consider the compact interval between them. The function \[
\tau\longmapsto
P_{\gamma,\tau\to t_0}^{\mathbf{F}_r}
\circ
\iota_{\dot\gamma(\tau)}
\] is continuous and takes values in the space of linear operators between finite-dimensional vector spaces. Thus there is a constant $C>0$ such that \[
\left\|
P_{\gamma,\tau\to t_0}^{\mathbf{F}_r}
\circ
\iota_{\dot\gamma(\tau)}
\right\|
\leq C
\] for every $\tau$ between $t_0$ and $t$.

Consequently, \[
\left|
\int_{t_0}^{t}
P_{\gamma,\tau\to t_0}^{\mathbf{F}_r}
\left(
\iota_{\dot\gamma(\tau)}
\bigl(
\nabla^{r+1}\boldsymbol{\phi}_n(\gamma(\tau))-\mathbf{A}_{r+1}(\gamma(\tau))
\bigr)
\right)
\,d\tau
\right|
\] \[
\qquad\leq
C\,|t-t_0|\,
\|\nabla^{r+1}\boldsymbol{\phi}_n-\mathbf{A}_{r+1}\|_{L^\infty(M,T^{(0,r+1)}(TM)\otimes\mathbf E)}.
\] The right-hand side converges to $0$, since $\nabla^{r+1}\boldsymbol{\phi}_n\to \mathbf{A}_{r+1}$ uniformly on $M$. Therefore \[
P_{\gamma,t\to t_0}^{\mathbf{F}_r}
\bigl(
\mathbf{A}_r(\gamma(t))
\bigr)
-
\mathbf{A}_r(\gamma(t_0))
\] \[
\qquad =
\int_{t_0}^{t}
P_{\gamma,\tau\to t_0}^{\mathbf{F}_r}
\left(
\iota_{\dot\gamma(\tau)}
\bigl(\mathbf{A}_{r+1}(\gamma(\tau))\bigr)
\right)
\,d\tau.
\]

This identity has exactly the form in Lemma~\ref{lema: criterio integral derivada covariante}, applied to the bundle $\mathbf{F}_r$, the section $\mathbf{A}_r$, and the tensor $\mathbf{A}_{r+1}$ viewed as an element of $\Gamma^0(T^*M\otimes \mathbf{F}_r)$. We conclude that $\mathbf{A}_r$ is of class $C^1$ and that $\nabla^{\mathbf{F}_r}\mathbf{A}_r=\mathbf{A}_{r+1}$.

Using this compatibility, we prove by induction on $r\in\{0,\dots,m\}$ that $\boldsymbol{\phi}\in\Gamma^r(\mathbf{E})$ and $\nabla^r\boldsymbol{\phi}=\mathbf{A}_r$.

For $r=0$, the assertion is immediate, since $\boldsymbol{\phi}=\mathbf{A}_0$ and $\mathbf{A}_0$ is continuous.

Now suppose the assertion holds for some $r\in\{0,\dots,m-1\}$. Then $\boldsymbol{\phi}\in\Gamma^r(\mathbf{E})$ and $\nabla^r\boldsymbol{\phi}=\mathbf{A}_r$. Since we have just proved that $\mathbf{A}_r$ is of class $C^1$ and $\nabla^{\mathbf{F}_r}\mathbf{A}_r=\mathbf{A}_{r+1}$, the definition of higher-order covariant derivative gives $\boldsymbol{\phi}\in\Gamma^{r+1}(\mathbf{E})$ and $\nabla^{r+1}\boldsymbol{\phi}=\mathbf{A}_{r+1}$.

By induction, $\boldsymbol{\phi}\in\Gamma^m(\mathbf{E})$ and $\nabla^r\boldsymbol{\phi}=\mathbf{A}_r$ for every $r\in\{0,\dots,m\}$.

Since we also have $\operatorname{supp}(\boldsymbol{\phi})\subseteq K$, we conclude that $\boldsymbol{\phi}\in\Gamma_K^m(\mathbf{E})$.

Finally, for each $r\in\{0,\dots,m\}$, $\nabla^r\boldsymbol{\phi}_n\to\nabla^r\boldsymbol{\phi}$ uniformly on $K$. Therefore \[
\|\boldsymbol{\phi}_n-\boldsymbol{\phi}\|_{C^m(K)}
=
\max_{0\leq r\leq m}
\|\nabla^r\boldsymbol{\phi}_n-\nabla^r\boldsymbol{\phi}\|_{C^0(K,T^{(0,r)}(TM)\otimes\mathbf E)}
\longrightarrow 0.
\] Thus every Cauchy sequence in $\Gamma_K^m(\mathbf{E})$ converges in $\Gamma_K^m(\mathbf{E})$. We conclude that $\bigl(\Gamma_K^m(\mathbf{E}),\|\cdot\|_{C^m(K)}\bigr)$ is a Banach space. \end{proof}

Having constructed the Banach space $\Gamma_K^m(\mathbf{E})$ for each finite order, we pass to the space of smooth sections supported in $K$. This space is described by the countable family of norms $C^m$, rather than by a single norm.

\begin{definition}\label{def:topologia-GammaK}\index{smooth sections!topology on a compact support} Let $(M,\mathbf{g})$ be a Riemannian manifold with or without boundary, let $\mathbf{E}\to M$ be a smooth vector bundle equipped with a bundle metric $\mathbf{h}_{\mathbf{E}}$, and let $K\subseteq M$ be compact. Equip $\Gamma_K(\mathbf{E})$ with the locally convex topology induced by the family of seminorms \[
\mathbf{u}\longmapsto \|\mathbf{u}\|_{C^m(K)},
\qquad m\in\mathbb N_0.
\] This is the topology of uniform convergence on $K$ of sections and all their covariant derivatives. \end{definition}

\begin{lemma}\label{lema: base local GammaK} Let $(M,\mathbf{g})$ be a Riemannian manifold with or without boundary, let $\mathbf{E}\to M$ be a smooth vector bundle equipped with a bundle metric $\mathbf{h}_{\mathbf{E}}$, and let $K\subseteq M$ be compact. Then the topology on $\Gamma_K(\mathbf{E})$ generated by the family of seminorms $(\|\cdot\|_{C^m(K)})_{m\in\mathbb N_0}$ has as a local base at $0$ the family of sets \[
V_{K,m,\ell}
:=
\left\{
\mathbf{u}\in \Gamma_K(\mathbf{E})\ \middle|\
\|\mathbf{u}\|_{C^m(K)}<\frac{1}{\ell}
\right\},
\qquad m\in\mathbb N_0,\quad \ell\in\mathbb N.
\] Consequently, a neighborhood base at $\mathbf{u}\in \Gamma_K(\mathbf{E})$ is given by \[
\mathbf{u}+V_{K,m,\ell}
=
\left\{
\mathbf{v}\in \Gamma_K(\mathbf{E})\ \middle|\
\|\mathbf{u}-\mathbf{v}\|_{C^m(K)}<\frac{1}{\ell}
\right\}.
\] \end{lemma}

\begin{proof} By Definition~\ref{def: topologia inducida por seminormas} of the topology induced by seminorms, a local base at $0$ consists of sets of the form \[
U
=
\left\{
\mathbf{u}\in \Gamma_K(\mathbf{E})\ \middle|\
\|\mathbf{u}\|_{C^{m_1}(K)}<\varepsilon_1,\dots,\|\mathbf{u}\|_{C^{m_q}(K)}<\varepsilon_q
\right\}.
\] Let $m:=\displaystyle\max\{m_1,\dots,m_q\}$. Then, since the norms $C^m$ increase with the order, \[
\|\mathbf{u}\|_{C^{m_i}(K)}\leq \|\mathbf{u}\|_{C^m(K)}
\qquad \forall i\in\{1,\dots,q\}.
\] Choose $\ell\in\mathbb N$ such that $\displaystyle\frac{1}{\ell}<\min\{\varepsilon_1,\dots,\varepsilon_q\}$. If $\mathbf{u}\in V_{K,m,\ell}$, then \[
\|\mathbf{u}\|_{C^{m_i}(K)}
\leq
\|\mathbf{u}\|_{C^m(K)}
<
\frac{1}{\ell}
<
\varepsilon_i
\] for every $i\in\{1,\dots,q\}$. Therefore $V_{K,m,\ell}\subseteq U$.

This proves that the family $V_{K,m,\ell}$ is a local base at $0$. The general case follows by translation, since translations are homeomorphisms of a topological vector space. \end{proof}

\begin{proposition}\label{prop: GammaK Frechet} Let $(M,\mathbf{g})$ be a Riemannian manifold with or without boundary, let $\mathbf{E}\to M$ be a smooth vector bundle equipped with a bundle metric $\mathbf{h}_{\mathbf{E}}$, and let $K\subseteq M$ be compact. For $\mathbf{u},\mathbf{v}\in\Gamma_K(\mathbf{E})$, define \[
d_K(\mathbf{u},\mathbf{v})
:=
\displaystyle\sum_{m=0}^{\infty}
\frac{1}{2^m}
\frac{\|\mathbf{u}-\mathbf{v}\|_{C^m(K)}}{1+\|\mathbf{u}-\mathbf{v}\|_{C^m(K)}}.
\] Then the following properties hold: \begin{enumerate}[label=(\alph*)] \item $d_K$ is a translation-invariant metric on $\Gamma_K(\mathbf{E})$. \item The topology on $\Gamma_K(\mathbf{E})$ generated by the family of seminorms $(\|\cdot\|_{C^m(K)})_{m\in\mathbb N_0}$ agrees with the topology induced by $d_K$. \item $(\Gamma_K(\mathbf{E}),d_K)$ is complete. \end{enumerate} In particular, $\Gamma_K(\mathbf{E})$ is a Fréchet space. \end{proposition}

\begin{proof} Consider the countable family of seminorms \[
p_m(\mathbf{u}):=\|\mathbf{u}\|_{C^m(K)},
\qquad m\in\mathbb N_0,
\] defined on $\Gamma_K(\mathbf{E})$. In fact, each $p_m$ is a norm on $\Gamma_K(\mathbf{E})$: if $p_m(\mathbf{u})=0$, then in particular $\mathbf{u}=0$ on $K$; since $\operatorname{supp}(\mathbf{u})\subseteq K$, it follows that $\mathbf{u}=0$ on $M\setminus K$, and hence $\mathbf{u}=0$ throughout $M$. Thus the family $(p_m)_{m\in\mathbb N_0}$ separates points.

Properties $(a)$ and $(b)$ are therefore a special case of Proposition~\ref{prop: metrizacion por seminormas}. The metric constructed there, applied to the family $(p_m)_{m\in\mathbb N_0}$, is precisely $d_K$; hence $d_K$ is translation invariant and induces the locally convex topology generated by the norms $\|\cdot\|_{C^m(K)}$.

It remains to prove completeness. We use sequential completeness of projective intersections. First observe that, within the Hausdorff topological vector space $Y:=\Gamma_K^0(\mathbf{E})$, we have the algebraic equality \[
\Gamma_K(\mathbf{E})
=
\bigcap_{m\in\mathbb N_0}\Gamma_K^m(\mathbf{E}).
\] If a section belongs to $\Gamma_K(\mathbf{E})$, it is smooth and therefore belongs to $\Gamma_K^m(\mathbf{E})$ for every $m\in\mathbb N_0$. Conversely, if a section belongs to $\Gamma_K^m(\mathbf{E})$ for every $m\in\mathbb N_0$, it is of class $C^m$ for every $m$, hence smooth, and its support is contained in $K$.

For each $m\in\mathbb N_0$, the space $\Gamma_K^m(\mathbf{E})$ is Banach by Proposition~\ref{prop: GammaKm es Banach}. Moreover, the inclusion \[
\Gamma_K^m(\mathbf{E})\hookrightarrow \Gamma_K^0(\mathbf{E})
\] is continuous, since \[
\|\mathbf{u}\|_{C^0(K)}\leq \|\mathbf{u}\|_{C^m(K)}.
\] On the other hand, the topology on $\Gamma_K(\mathbf{E})$ generated by the norms $\|\cdot\|_{C^m(K)}$ agrees with the projective topology associated with the natural inclusions \[
\iota_m\colon \Gamma_K(\mathbf{E})\longrightarrow \Gamma_K^m(\mathbf{E}),
\qquad m\in\mathbb N_0.
\] This is because that topology is, by definition, the initial topology making all the inclusions $\iota_m$ continuous, or equivalently the topology generated by the seminorms $\mathbf{u}\mapsto\|\iota_m(\mathbf{u})\|_{C^m(K)}=\|\mathbf{u}\|_{C^m(K)}$; see Definition~\ref{def: topologia inicial} and Definition~\ref{def: topologia inducida por seminormas}.

Applying Proposition~\ref{prop: completitud secuencial limite proyectivo} and reindexing the family $\{\Gamma_K^m(\mathbf{E})\}_{m\in\mathbb N_0}$ by $n=m+1$, we conclude that $\Gamma_K(\mathbf{E})$ is sequentially complete in this projective topology. Since this same topology is induced by the translation-invariant metric $d_K$, a sequence is Cauchy with respect to $d_K$ if and only if it is Cauchy in the corresponding vector topology. By sequential completeness, every Cauchy sequence for $d_K$ converges in $\Gamma_K(\mathbf{E})$; since $d_K$ induces the topology, this is convergence in the metric $d_K$. Therefore $(\Gamma_K(\mathbf{E}),d_K)$ is complete.

Finally, $\Gamma_K(\mathbf{E})$ is locally convex because its topology is generated by seminorms, Hausdorff because the family $(p_m)_{m\in\mathbb N_0}$ separates points, metrizable by $d_K$, and complete by what we have just proved. By the definition of a Fréchet space, or equivalently by Theorem~\ref{teo: caracterizacion Frechet por seminormas}, we conclude that $\Gamma_K(\mathbf{E})$ is a Fréchet space. \end{proof}

\begin{proposition}\label{prop:compatibilidad-soportes-compactos-secciones} Let $(M,\mathbf{g})$ be a Riemannian manifold with or without boundary, let $\mathbf{E}\to M$ be a smooth vector bundle equipped with a bundle metric $\mathbf{h}_{\mathbf{E}}$, and let $K,L\subseteq M$ be compact sets such that $K\subseteq L$. Then: \begin{enumerate}[label=(\alph*)] \item for each $m\in\mathbb N_0$, the natural inclusion $\Gamma_K^m(\mathbf{E})\hookrightarrow\Gamma_L^m(\mathbf{E})$ is an isometry onto its image, that is, \[
 \|\mathbf{u}\|_{C^m(L)}=\|\mathbf{u}\|_{C^m(K)}
 \qquad
 \forall \mathbf{u}\in\Gamma_K^m(\mathbf{E});
 \] \item the topology on $\Gamma_K^m(\mathbf{E})$ agrees with the subspace topology induced by $\Gamma_L^m(\mathbf{E})$; \item the Fréchet topology on $\Gamma_K(\mathbf{E})$ agrees with the subspace topology induced by $\Gamma_L(\mathbf{E})$; \item the images of $\Gamma_K^m(\mathbf{E})$ in $\Gamma_L^m(\mathbf{E})$ and of $\Gamma_K(\mathbf{E})$ in $\Gamma_L(\mathbf{E})$ are closed. \end{enumerate} \end{proposition}

\begin{proof} Let $\mathbf{u}\in\Gamma_K^m(\mathbf{E})$. Since $\operatorname{supp}(\mathbf{u})\subseteq K$, we have $\mathbf{u}=0$ on the open subset $M\setminus K$. Thus, for each $j\in\{0,\dots,m\}$, we also have $\nabla^j \mathbf{u}=0$ on $M\setminus K$. In particular, on $L\setminus K$ all covariant derivatives of $\mathbf{u}$ up to order $m$ vanish. This gives \[
\sup_{x\in L}|\nabla^j \mathbf{u}(x)|_{\mathbf{g},\mathbf{h}_{\mathbf{E}}}
=
\sup_{x\in K}|\nabla^j \mathbf{u}(x)|_{\mathbf{g},\mathbf{h}_{\mathbf{E}}}
\qquad
0\leq j\leq m.
\] Taking the maximum over $0\leq j\leq m$ yields $\|\mathbf{u}\|_{C^m(L)}=\|\mathbf{u}\|_{C^m(K)}$. This proves $(a)$.

Part $(b)$ follows immediately from $(a)$, since the norm on $\Gamma_L^m(\mathbf{E})$ restricted to the subspace $\Gamma_K^m(\mathbf{E})$ agrees with the intrinsic norm of $\Gamma_K^m(\mathbf{E})$.

For $(c)$, recall that the topology on $\Gamma_K(\mathbf{E})$ is induced by the family of norms $\|\cdot\|_{C^m(K)}$, with $m\in\mathbb N_0$, whereas the topology on $\Gamma_L(\mathbf{E})$ is induced by the family $\|\cdot\|_{C^m(L)}$, with $m\in\mathbb N_0$. By $(a)$, restricting the latter norms to $\Gamma_K(\mathbf{E})$ gives exactly the family $\|\cdot\|_{C^m(K)}$. Therefore the topologies agree.

For $(d)$, let $\mathbf{u}_j\in\Gamma_K^m(\mathbf{E})$ be a sequence converging to $\mathbf{u}\in\Gamma_L^m(\mathbf{E})$. In particular, $\mathbf{u}_j\to \mathbf{u}$ uniformly on $L$. Since $\mathbf{u}_j=0$ on $L\setminus K$, we also have $\mathbf{u}=0$ there; outside $L$ the section $\mathbf{u}$ already vanishes, so $\operatorname{supp}\mathbf{u}\subseteq K$. This proves closedness of the image at the $C^m$ level. If $\mathbf{u}_j\in\Gamma_K(\mathbf{E})$ converges to $\mathbf{u}$ in $\Gamma_L(\mathbf{E})$, the seminorm $C^0(L)$ gives the same argument. Thus both images are closed. \end{proof}

\begin{definition}\label{def: topologias secciones}\index{test sections!topology} Let $M$ be a smooth manifold with or without boundary, and let $\mathbf{E}\longrightarrow M$ be a smooth vector bundle. Fix an increasing exhaustion of $M$ by compact sets \[
K_1\subseteq K_2\subseteq \cdots \subseteq M,
\qquad
K_n\subseteq \operatorname{Int}(K_{n+1}),
\qquad
M=\bigcup_{n=1}^{\infty}K_n.
\]

\begin{enumerate}[label=(\alph*)] \item For each $m\in\mathbb N_0$, equip the space $\Gamma^m(\mathbf{E})$ with the locally convex topology generated by the family of seminorms \[
 p_{n,j}(u)
 :=
 \sup_{x\in K_n}|\nabla^ju(x)|_{\mathbf{g},\mathbf{h}_{\mathbf{E}}},
 \qquad n\in\mathbb N,\quad 0\leq j\leq m.
 \]

\item Equip the space $\Gamma(\mathbf{E})$ with the locally convex topology generated by the family of seminorms \[
 p_{n,m}(u)
 :=
 \max_{0\leq j\leq m}\sup_{x\in K_n}|\nabla^ju(x)|_{\mathbf{g},\mathbf{h}_{\mathbf{E}}},
 \qquad n\in\mathbb N,\quad m\in\mathbb N_0.
 \]

\item For each $m\in\mathbb N_0$, equip the space $\Gamma_c^m(\mathbf{E})$ with the final locally convex topology associated with the inclusions \[
 \Gamma_{K_n}^m(\mathbf{E})\hookrightarrow\Gamma_c^m(\mathbf{E}),
 \qquad n\in\mathbb N.
 \] In other words, \[
 \Gamma_c^m(\mathbf{E})
 =
 \varinjlim_{n\in\mathbb N}\Gamma_{K_n}^m(\mathbf{E}).
 \] By Proposition~\ref{prop:compatibilidad-soportes-compactos-secciones}, the inclusions are embeddings with closed image. Since each $\Gamma_{K_n}^m(\mathbf{E})$ is a Banach space, $\Gamma_c^m(\mathbf{E})$ is a strict LB inductive limit and, in particular, a strict LF space.

\item Equip the space $\Gamma_c(\mathbf{E})$ with the final locally convex topology associated with the inclusions \[
 \Gamma_{K_n}(\mathbf{E})\hookrightarrow\Gamma_c(\mathbf{E}),
 \qquad n\in\mathbb N.
 \] That is, \[
 \Gamma_c(\mathbf{E})
 =
 \varinjlim_{n\in\mathbb N}\Gamma_{K_n}(\mathbf{E}).
 \] By Proposition~\ref{prop:compatibilidad-soportes-compactos-secciones}, the inclusions are embeddings with closed image. Since each $\Gamma_{K_n}(\mathbf{E})$ is a Fréchet space, $\Gamma_c(\mathbf{E})$ is a strict LF space. \end{enumerate} \end{definition}

\begin{proposition}\label{prop:independencia-topologias-secciones} Let $(M,\mathbf{g})$ be a Riemannian manifold with or without boundary, and let $\mathbf{E}\to M$ be a smooth vector bundle of finite rank. The topologies defined in \ref{def: topologias secciones} are independent of the chosen metrics, connections, and increasing exhaustion by compact sets. \end{proposition}

\begin{proof} Independence of the metrics and auxiliary connections follows from Theorem~\ref{norma ck independiente de metrica y conexion}. On each fixed compact set $K\subseteq M$, the norms $C^m(K)$ constructed using two different choices of metrics and connections are equivalent. Thus, for each fixed $m$, both choices generate the same topology on $\Gamma_K^m(\mathbf{E})$ and $\Gamma_K(\mathbf{E})$.

It remains to justify independence of the exhaustion. Let $(K_n)_{n\in\mathbb N}$ and $(L_q)_{q\in\mathbb N}$ be two increasing exhaustions of $M$ by compact sets, with $K_n\subseteq \operatorname{Int}(K_{n+1})$, $L_q\subseteq \operatorname{Int}(L_{q+1})$, $\displaystyle M=\bigcup_{n=1}^{\infty}K_n$, and $\displaystyle M=\bigcup_{q=1}^{\infty}L_q$.

First show that each compact set in one exhaustion lies in a compact set of the other. Fix $n\in\mathbb N$. Since $\displaystyle M=\bigcup_{q=1}^{\infty}L_q$ and $L_q\subseteq \operatorname{Int}(L_{q+1})$, the open subsets $\operatorname{Int}(L_q)$ cover $M$. In particular, they cover the compact set $K_n$. By compactness, there exist $q_1,\dots,q_s$ such that \[
K_n\subseteq \bigcup_{\alpha=1}^{s}\operatorname{Int}(L_{q_\alpha}).
\] If $q_0:=\displaystyle\max\{q_1,\dots,q_s\}$, then, since the family $(L_q)$ is increasing, \[
K_n\subseteq L_{q_0}.
\] The same argument, interchanging $K$ and $L$, shows that for each $q\in\mathbb N$ there exists $n_0\in\mathbb N$ such that $L_q\subseteq K_{n_0}$.

Now consider the topology on $\Gamma^m(\mathbf{E})$. For the exhaustion $(K_n)$, it is generated by the seminorms \[
p_{n,j}^{K}(\mathbf{u}):=\sup_{x\in K_n}|\nabla^j \mathbf{u}(x)|_{\mathbf{g},\mathbf{h}_{\mathbf{E}}},
\qquad n\in\mathbb N,\quad 0\leq j\leq m,
\] whereas for the exhaustion $(L_q)$, it is generated by the seminorms \[
p_{q,j}^{L}(\mathbf{u}):=\sup_{x\in L_q}|\nabla^j \mathbf{u}(x)|_{\mathbf{g},\mathbf{h}_{\mathbf{E}}},
\qquad q\in\mathbb N,\quad 0\leq j\leq m.
\] If $K_n\subseteq L_{q_0}$, then $p_{n,j}^{K}(\mathbf{u})\leq p_{q_0,j}^{L}(\mathbf{u})$ for every $\mathbf{u}\in\Gamma^m(\mathbf{E})$. Thus every basic neighborhood defined using the exhaustion $(K_n)$ contains a basic neighborhood defined using the exhaustion $(L_q)$. The reverse argument, using the fact that each $L_q$ lies in some $K_{n_0}$, gives the reverse inclusion between the topologies. Thus both exhaustions induce the same topology on $\Gamma^m(\mathbf{E})$. The same reasoning applies to $\Gamma(\mathbf{E})$, using the seminorms \[
p_{n,r}^{K}(\mathbf{u}):=\max_{0\leq j\leq r}\sup_{x\in K_n}|\nabla^j \mathbf{u}(x)|_{\mathbf{g},\mathbf{h}_{\mathbf{E}}}
\] and their counterparts associated with $(L_q)$.

We now prove independence of the exhaustion for the spaces with compact support. Denote by $\tau_K^m$ the final locally convex topology on $\Gamma_c^m(\mathbf{E})$ associated with the family $(\Gamma_{K_n}^m(\mathbf{E}))_{n\in\mathbb N}$, and by $\tau_L^m$ the final locally convex topology associated with the family $(\Gamma_{L_q}^m(\mathbf{E}))_{q\in\mathbb N}$.

We prove $\tau_K^m=\tau_L^m$ using the maximality property in Theorem~\ref{teo: topologia-final-lc}. For each $n\in\mathbb N$, choose $q_0\in\mathbb N$ such that $K_n\subseteq L_{q_0}$. Then the inclusion \[
\Gamma_{K_n}^m(\mathbf{E})\hookrightarrow \Gamma_c^m(\mathbf{E})
\] viewed with the topology $\tau_L^m$ is continuous, since it factors as \[
\Gamma_{K_n}^m(\mathbf{E})
\hookrightarrow
\Gamma_{L_{q_0}}^m(\mathbf{E})
\hookrightarrow
(\Gamma_c^m(\mathbf{E}),\tau_L^m),
\] where the first arrow is continuous by Proposition~\ref{prop:compatibilidad-soportes-compactos-secciones} and the second by the definition of $\tau_L^m$ as a final locally convex topology.

Thus $\tau_L^m$ is a locally convex topology on $\Gamma_c^m(\mathbf{E})$ making all the inclusions $\Gamma_{K_n}^m(\mathbf{E})\hookrightarrow\Gamma_c^m(\mathbf{E})$ continuous. Since $\tau_K^m$ is the finest topology with this property, we obtain $\tau_L^m\subseteq\tau_K^m$. Interchanging the exhaustions $(K_n)$ and $(L_q)$ gives $\tau_K^m\subseteq\tau_L^m$. Consequently, $\tau_K^m=\tau_L^m$.

The same argument proves independence of the topology on $\Gamma_c(\mathbf{E})$. If $K_n\subseteq L_{q_0}$, the inclusion \[
\Gamma_{K_n}(\mathbf{E})\hookrightarrow\Gamma_{L_{q_0}}(\mathbf{E})
\] is continuous by Proposition~\ref{prop:compatibilidad-soportes-compactos-secciones}; hence the final locally convex topology defined by the exhaustion $(L_q)$ makes all the inclusions $\Gamma_{K_n}(\mathbf{E})\hookrightarrow\Gamma_c(\mathbf{E})$ continuous. By maximality in Theorem~\ref{teo: topologia-final-lc}, one topology contains the other. Interchanging the exhaustions gives the reverse containment, so the two agree. \end{proof}

\begin{definition}\label{def:espacio-DME}\index{test sections!of a vector bundle} Let $M$ be a smooth manifold with or without boundary, and let $\mathbf{E}\longrightarrow M$ be a smooth vector bundle. Define \[
\mathcal D(M,\mathbf{E}):=\Gamma_c(\mathbf{E}),
\] and equip this space with the strict locally convex inductive topology \[
\mathcal D(M,\mathbf{E})
=
\varinjlim_{n\in\mathbb N}\Gamma_{K_n}(\mathbf{E}),
\] where $(K_n)_{n\in\mathbb N}$ is an increasing exhaustion of $M$ by compact sets. By Proposition~\ref{prop:independencia-topologias-secciones}, this topology is independent of the chosen exhaustion. In particular, $\mathcal D(M,\mathbf{E})$ is an LF space. \end{definition}

The LF structure of $\mathcal D(M,\mathbf{E})$ gives a concrete description of sequential convergence: a sequence converges if and only if all its terms are supported in a common compact set and converge uniformly there together with all their covariant derivatives.

\begin{corollary}\label{cor:convergencia-en-DME} Let $M$ be a smooth manifold with or without boundary, and let $\mathbf{E}\to M$ be a smooth vector bundle. Let $\boldsymbol{\phi}_n,\boldsymbol{\phi}\in\mathcal D(M,\mathbf{E})$, with $n\in\mathbb N$. Then $\boldsymbol{\phi}_n\to\boldsymbol{\phi}$ in $\mathcal D(M,\mathbf{E})$ if and only if there exists a compact set $K\subseteq M$ such that: \begin{enumerate}[label=(\alph*)] \item $\operatorname{supp}(\boldsymbol{\phi}_n)\subseteq K$ for every $n\in\mathbb N$ and $\operatorname{supp}(\boldsymbol{\phi})\subseteq K$; \item for every $k\in\mathbb N_0$, \[
 \|\boldsymbol{\phi}_n-\boldsymbol{\phi}\|_{C^k(K)}\to0
 \qquad \text{when } n\to\infty.
 \] \end{enumerate} \end{corollary}

\begin{proof} By definition, $\mathcal D(M,\mathbf{E})$ is the LF space generated by the sequence of Fréchet spaces $\Gamma_{K_n}(\mathbf{E})$. By the corresponding part of Theorem~\ref{teo: propiedades LF}, the sequence $\boldsymbol{\phi}_n$ converges to $\boldsymbol{\phi}$ in $\mathcal D(M,\mathbf{E})$ if and only if there exists an index $N\in\mathbb N$ such that \[
\{\boldsymbol{\phi}\}\cup\{\boldsymbol{\phi}_n\mid n\in\mathbb N\}\subseteq \Gamma_{K_N}(\mathbf{E})
\] and, moreover, $\boldsymbol{\phi}_n\to\boldsymbol{\phi}$ in the Fréchet topology of $\Gamma_{K_N}(\mathbf{E})$.

The first condition means that all supports lie in the compact set $K_N$. Taking $K=K_N$ gives $(a)$. On the other hand, the topology on $\Gamma_{K_N}(\mathbf{E})$ is induced by the family of seminorms $\|\cdot\|_{C^k(K_N)}$, with $k\in\mathbb N_0$. By the characterization of convergence in seminorm-induced topologies, Proposition~\ref{prop: convergencia y Cauchy por seminormas}, convergence $\boldsymbol{\phi}_n\to\boldsymbol{\phi}$ in $\Gamma_{K_N}(\mathbf{E})$ is equivalent to \[
\|\boldsymbol{\phi}_n-\boldsymbol{\phi}\|_{C^k(K_N)}\to0
\qquad \text{for every }k\in\mathbb N_0.
\] This proves $(b)$.

Conversely, suppose there exists a compact set $K\subseteq M$ satisfying $(a)$ and $(b)$. Since $(K_n)_{n\in\mathbb N}$ is an increasing exhaustion with $K_n\subseteq\operatorname{Int}(K_{n+1})$, there exists $N\in\mathbb N$ such that $K\subseteq K_N$. Then $\boldsymbol{\phi}_n,\boldsymbol{\phi}\in\Gamma_{K_N}(\mathbf{E})$ for every $n$.

Fix $k\in\mathbb N_0$. Since the sections $\boldsymbol{\phi}_n-\boldsymbol{\phi}$ are supported in $K$, they vanish on $M\setminus K$. Thus all their covariant derivatives also vanish on $M\setminus K$. In particular, \[
\|\boldsymbol{\phi}_n-\boldsymbol{\phi}\|_{C^k(K_N)}
=
\|\boldsymbol{\phi}_n-\boldsymbol{\phi}\|_{C^k(K)}.
\] By hypothesis $(b)$, the right-hand side tends to $0$. Since this holds for every $k\in\mathbb N_0$, we obtain $\boldsymbol{\phi}_n\to\boldsymbol{\phi}$ in $\Gamma_{K_N}(\mathbf{E})$. Applying Theorem~\ref{teo: propiedades LF} once again, we conclude that $\boldsymbol{\phi}_n\to\boldsymbol{\phi}$ in $\mathcal D(M,\mathbf{E})$. \end{proof}

\section{Distributions on vector bundles} With these tools, we are ready to define distributions in the setting of vector bundles. \begin{definition}\label{def:distribuciones-en-haces-distribuciones}\index{distribution!on a vector bundle}\index{vector bundle!distribution}\glsadd{distribucion-haz} Let $M$ be a smooth manifold with or without boundary, let $\mathbf{E}\to M$ be a smooth vector bundle, and let $W$ be a finite-dimensional $\mathbb{K}$-vector space $(\mathbb{K}\in\{\mathbb{R},\mathbb{C}\})$.

Denote by $\mathcal{D}'(M,\mathbf{E},W)$ the space of all continuous $\mathbb{K}$-linear maps \[
\mathbf{F}\colon \mathcal{D}(M,\mathbf{E}^{*})\longrightarrow W
\].

Equip $\mathcal{D}'(M,\mathbf{E},W)$ with the weak-$*$ topology, defined as the initial topology associated with the family of evaluation maps \[
\operatorname{ev}_{\boldsymbol{\phi}}\colon \mathcal{D}'(M,\mathbf{E},W)\longrightarrow W,\qquad
\operatorname{ev}_{\boldsymbol{\phi}}(\mathbf{F})=\mathbf{F}(\boldsymbol{\phi}),
\qquad \boldsymbol{\phi}\in \mathcal{D}(M,\mathbf{E}^{*}),
\] that is, the topology generated by sets of the form \[
\operatorname{ev}_{\boldsymbol{\phi}}^{-1}(U)
=
\{\mathbf{F}\in \mathcal{D}'(M,\mathbf{E},W)\mid \mathbf{F}(\boldsymbol{\phi})\in U\},
\] where $\boldsymbol{\phi}\in \mathcal{D}(M,\mathbf{E}^{*})$ and $U\subseteq W$ is open.

Elements of $\mathcal{D}'(M,\mathbf{E},W)$ are called \textbf{distributions} on $\mathbf{E}$ with values in $W$.

When $W=\mathbb{K}$, we simply write $\mathcal{D}'(M,\mathbf{E})$, and if in addition $\mathbf{E}=T^{(0,0)}(TM)$, we write $\mathcal{D}'(M)$. \end{definition} \begin{remark}\label{obs: top debil * distribuciones haces} The weak-$*$ topology on $\mathcal{D}'(M,\mathbf{E},W)$ agrees with the topology of pointwise convergence on $\mathcal{D}(M,\mathbf{E}^{*})$. Thus a sequence $(\mathbf{F}_n)_{n\in\mathbb{N}}\subseteq \mathcal{D}'(M,\mathbf{E},W)$ converges to $\mathbf{F}\in \mathcal{D}'(M,\mathbf{E},W)$ if and only if \[
\mathbf{F}_n(\boldsymbol{\phi})\to \mathbf{F}(\boldsymbol{\phi})
\qquad \forall \boldsymbol{\phi}\in \mathcal{D}(M,\mathbf{E}^{*}).
\] More generally, a net $(\mathbf{F}_\alpha)_\alpha$ converges to $\mathbf{F}$ if and only if \[
\mathbf{F}_\alpha(\boldsymbol{\phi})\to \mathbf{F}(\boldsymbol{\phi})
\qquad \forall \boldsymbol{\phi}\in \mathcal{D}(M,\mathbf{E}^{*}).
\] \end{remark} \begin{proposition}\label{prop:distribuciones-en-haces-espacio-vectorial-dimension-finita-isomorfismo-natural} Let $M$ be a smooth manifold with or without boundary, let $\mathbf{E}\to M$ be a smooth vector bundle, let $W$ be a finite-dimensional $\mathbb{K}$-vector space, and let $N:=\dim W$. Then there exists a natural isomorphism of topological vector spaces \[
\mathcal{D}'(M,\mathbf{E},W)\cong \mathcal{D}'(M,\mathbf{E})^N,
\] where the right-hand side is equipped with the product topology. \end{proposition}

\begin{proof} Let $\{w_1,\dots,w_N\}$ be a basis of $W$ and $\{w^1,\dots,w^N\}$ the induced dual basis. Given $\mathbf{F}\in\mathcal{D}'(M,\mathbf{E},W)$, define \[
\mathbf{F}_j:=w^j\circ \mathbf{F}\in\mathcal{D}'(M,\mathbf{E}), \qquad j\in\{1,\dots,N\}.
\] Then, for every $\boldsymbol{\phi}\in\mathcal{D}(M,\mathbf{E}^*)$, we have \[
\mathbf{F}(\boldsymbol{\phi})=\displaystyle\sum_{j=1}^N \mathbf{F}_j(\boldsymbol{\phi})\,w_j.
\] This defines a linear isomorphism \[
\mathcal{D}'(M,\mathbf{E},W)\longrightarrow \mathcal{D}'(M,\mathbf{E})^N.
\]

Moreover, by the definition of the weak-$*$ topology, convergence in $\mathcal{D}'(M,\mathbf{E},W)$ is equivalent to componentwise convergence, showing that this isomorphism is a homeomorphism. \end{proof} \begin{remark}\label{obs:distribuciones-en-haces-particular-identificacion-natural-cual-elegir-base} In particular, there is a natural identification \[
\mathcal{D}'(M,\mathbf{E},W)\cong \mathcal{D}'(M,\mathbf{E})\otimes W,
\] which, upon choosing a basis of $W$, reduces to the preceding identification with a finite product. \end{remark} Continuity can be characterized as follows, analogously to the Euclidean case. \begin{theorem}\label{teo: equivalencias de continuidad en D'(M,E,W)}\index{continuity criteria in d m e w@continuity criteria in D'(M,E,W)} Let $M$ be a smooth manifold with or without boundary, $\mathbf{E} \longrightarrow M$ a smooth vector bundle, and $W$ a finite-dimensional $\mathbb{K}$-vector space equipped with a norm $|\cdot|_W$. Let $\mathbf{F}\colon \mathcal{D}(M,\mathbf{E}^*)\longrightarrow W$ be a linear operator. Then the following properties are equivalent: \begin{enumerate}[label=(\alph*)] \item $\mathbf{F}$ is continuous. \item If $\{\boldsymbol{\phi}_{k}\}_{k\in \mathbb{N}}$ converges to $\boldsymbol{\phi}$ in $\mathcal{D}(M,\mathbf{E}^*)$, then \[
 \lim_{k\to \infty} \mathbf{F}(\boldsymbol{\phi}_{k}) = \mathbf{F}(\boldsymbol{\phi}).
 \] \item For every compact set $K\subseteq M$, the restriction $\mathbf{F}\restriction_{\mathcal{D}_{K}(M, \mathbf{E}^*)}$ is continuous. \item For every compact set $K\subseteq M$, there exist $m \in \mathbb{N}_{0}$ and $c_{K}>0$ such that \[
 |\mathbf{F}(\boldsymbol{\phi})|_W \leq c_{K}\|\boldsymbol{\phi}\|_{C^{m}(K)}\qquad \forall \boldsymbol{\phi}\in \mathcal{D}_{K}(M, \mathbf{E}^*).
 \] \end{enumerate} \end{theorem}

\begin{proof} (a)$\Rightarrow$(b). Let $\{\boldsymbol{\phi}_k\}_{k\in\mathbb{N}}$ be a sequence in $\mathcal{D}(M,\mathbf{E}^*)$ such that $\boldsymbol{\phi}_k\to \boldsymbol{\phi}$. Continuity of $\mathbf{F}$ directly gives \[
\mathbf{F}(\boldsymbol{\phi}_k)\longrightarrow \mathbf{F}(\boldsymbol{\phi}).
\]

\medskip

(b)$\Rightarrow$(c). Let $K\subseteq M$ be a compact set. Consider the canonical inclusion \[
\iota\colon \mathcal{D}_K(M,\mathbf{E}^*)\hookrightarrow \mathcal{D}(M,\mathbf{E}^*),
\] which is continuous by construction of the inductive topology.

Let $\{\boldsymbol{\phi}_k\}_{k\in\mathbb{N}}$ be a sequence in $\mathcal{D}_K(M,\mathbf{E}^*)$ such that $\boldsymbol{\phi}_k\to \boldsymbol{\phi}$ in $\mathcal{D}_K(M,\mathbf{E}^*)$. Then $\boldsymbol{\phi}_k\to \boldsymbol{\phi}$ in $\mathcal{D}(M,\mathbf{E}^*)$, and by (b), \[
\mathbf{F}(\boldsymbol{\phi}_k)\longrightarrow \mathbf{F}(\boldsymbol{\phi}).
\] This shows that the restriction $\mathbf{F}\restriction_{\mathcal{D}_K(M,\mathbf{E}^*)}$ is sequentially continuous. Since $\mathcal{D}_K(M,\mathbf{E}^*)$ is a Fréchet space, it is in particular metrizable, and Proposition~\ref{prop: continuidad secuencial en espacios metrizables} shows that sequential continuity implies continuity. Therefore $\mathbf{F}\restriction_{\mathcal{D}_K(M,\mathbf{E}^*)} \text{ is continuous.}$

(c)$\Rightarrow$(a). Let $\{K_n\}_{n\in\mathbb{N}}$ be an exhausting sequence of compact sets such that \[
\mathcal{D}(M,\mathbf{E}^*)=\varinjlim_{n\in\mathbb{N}} \mathcal{D}_{K_n}(M,\mathbf{E}^*).
\] By hypothesis, for every $n\in\mathbb{N}$ the restriction $\mathbf{F}\restriction_{\mathcal{D}_{K_n}(M,\mathbf{E}^*)}$ is continuous.

By the definition of the inductive topology (Theorem~\ref{teo: topologia-final-lc}), a linear operator \[
\mathbf{F}\colon \mathcal{D}(M,\mathbf{E}^*)\longrightarrow W
\] is continuous if and only if all its restrictions to the spaces $\mathcal{D}_{K_n}(M,\mathbf{E}^*)$ are continuous. We conclude that $\mathbf{F}$ is continuous.

\medskip

(c)$\Rightarrow$(d). Let $K\subseteq M$ be a compact set. Then $\mathcal{D}_K(M,\mathbf{E}^*)$ is a locally convex topological vector space whose topology is generated by the increasing family of seminorms \[
\|\boldsymbol{\phi}\|_{C^m(K)},\qquad m\in\mathbb{N}_0.
\] By (c), the restriction $\mathbf{F}\restriction_{\mathcal{D}_K(M,\mathbf{E}^*)}$ is continuous. Applying Theorem~\ref{teo: continuidad mediante seminormas} gives $m\in\mathbb{N}_0$ and $c_K>0$ such that \[
|\mathbf{F}(\boldsymbol{\phi})|_W \le c_K \|\boldsymbol{\phi}\|_{C^m(K)} \qquad \forall \boldsymbol{\phi}\in \mathcal{D}_K(M,\mathbf{E}^*).
\]

\medskip

(d)$\Rightarrow$(c). Let $K\subseteq M$ be compact. Suppose the estimate in (d) holds. If $\{\boldsymbol{\phi}_k\}_{k\in\mathbb{N}}$ converges to $\boldsymbol{\phi}$ in $\mathcal{D}_K(M,\mathbf{E}^*)$, then in particular \[
\|\boldsymbol{\phi}_k-\boldsymbol{\phi}\|_{C^m(K)}\longrightarrow 0
\] for every $m\in\mathbb{N}_0$. Applying the estimate in (d), we obtain \[
|\mathbf{F}(\boldsymbol{\phi}_k)-\mathbf{F}(\boldsymbol{\phi})|_W
\le c_K \|\boldsymbol{\phi}_k-\boldsymbol{\phi}\|_{C^m(K)}\longrightarrow 0.
\] This shows that $\mathbf{F}\restriction_{\mathcal{D}_K(M,\mathbf{E}^*)}$ is continuous.

\end{proof} With this in mind, we define the order of a distribution:

\begin{definition}\label{def: orden distribucion}\index{order of a distribution on a bundle@order of a distribution on a bundle} Let $M$ be a smooth manifold with or without boundary, let $\mathbf{E}\to M$ be a smooth vector bundle, let $W$ be a finite-dimensional vector space, and let $\mathbf{F}\colon \mathcal{D}(M,\mathbf{E}^*)\longrightarrow W$ be a distribution.

Let $K\subseteq M$ be a compact set. By part (d) of Theorem~\ref{teo: equivalencias de continuidad en D'(M,E,W)}, there exist $m_K\in \mathbb{N}_{0}$ and $C_K>0$ such that \[
|\mathbf{F}(\boldsymbol{\phi})|_W \leq C_K\|\boldsymbol{\phi}\|_{C^{m_K}(K)} \qquad \forall \boldsymbol{\phi}\in \mathcal{D}_{K}(M, \mathbf{E}^*).
\] The least integer $m_K$ satisfying this inequality is called the \textit{order} of $\mathbf{F}$ on the compact set $K$.

We say that $\mathbf{F}$ is a distribution of \textit{order} $m$ on $M$ if $m$ is the least nonnegative integer such that, for every compact set $K\subseteq M$, there exists a constant $C_K>0$ with \[
|\mathbf{F}(\boldsymbol{\phi})|_W \leq C_K\|\boldsymbol{\phi}\|_{C^{m}(K)} \qquad \forall \boldsymbol{\phi}\in \mathcal{D}_{K}(M, \mathbf{E}^*).
\] If no such integer $m$ exists, we say that $\mathbf{F}$ is a distribution of infinite order. \end{definition}

\begin{lemma}[Approximation with controlled enlargement of the support] \label{lem:densidad-soporte-aumentado-Cm} Let $M$ be a smooth manifold with or without boundary, and let $\mathbf{E}\to M$ be a smooth vector bundle. If $K\subseteq\operatorname{Int}L$ are compact, then for each $\mathbf{u}\in\Gamma_K^m(\mathbf{E})$ there exists a sequence $\mathbf{u}_j\in\Gamma_L(\mathbf{E})$ such that $\mathbf{u}_j\to \mathbf{u}$ in the $C^m(L)$ norm. \end{lemma}

\begin{proof} If $\dim M=0$, every section is smooth, and it suffices to take $\mathbf{u}_j=\mathbf{u}$. Suppose that $\dim M=n\geq1$. Choose an open subset $V$ with $K\subseteq V\Subset\operatorname{Int}L$, a finite cover of $\overline V$ by trivializing charts relatively compact in $\operatorname{Int}L$, and functions $\rho_a\in C_c^\infty(M)$ subordinate to these charts, indexed by $a\in\{1,\dots,N\}$, whose sum is one on a neighborhood of $K$. Then $\mathbf{u}=\displaystyle\sum_{a=1}^{N}\rho_a\mathbf{u}$. In each chart, the support of the components of $\rho_a\mathbf{u}$ is a compact subset of its coordinate domain; in particular, it is separated from its artificial boundaries.

Fix a nonnegative mollifier $\eta\in C_c^\infty(B(0,1))$ with integral one, and set $\eta_\varepsilon(z)=\varepsilon^{-n}\eta(\displaystyle\frac{z}{\varepsilon})$. In an interior chart, each component $f$ extends by zero to $\mathbb R^n$ as a compactly supported function of class $C^m$. The convolution $\eta_\varepsilon*f$ converges to $f$ uniformly together with derivatives of order at most $m$.

In a boundary chart, first extend $f$ by zero within $\mathbb H^n=\{x\in\mathbb R^n\mid x_n\geq0\}$, only across the artificial boundaries of the chart. This extension remains of class $C^m$ up to the boundary and has compact support in $\mathbb H^n$. Denote by $f_0$ its measurable extension by zero to the lower half-space, and define \[
f_\varepsilon(x):=(\eta_\varepsilon*f_0)(x+2\varepsilon e_n),
\qquad x\in\mathbb H^n.
\] The convolution is smooth on $\mathbb R^n$. Moreover, for $x_n\geq0$ and $z\in\operatorname{supp}\eta$, the last coordinate of $x+2\varepsilon e_n-\varepsilon z$ is at least $\varepsilon$. Thus, for every $\alpha\in\mathbb N_0^n$ with $|\alpha|\leq m$, we may differentiate using only interior values of $f$, obtaining \[
D^\alpha f_\varepsilon(x)
=
\int_{\mathbb R^n}\eta(z)
D^\alpha f(x+2\varepsilon e_n-\varepsilon z)\,dz.
\] The derivatives $D^\alpha f$ are uniformly continuous on the half-space, since they are continuous up to the boundary and compactly supported. Since $|2\varepsilon e_n-\varepsilon z|\leq3\varepsilon$ in this integral and $\int\eta=1$, it follows that \[
\max_{|\alpha|\leq m}\sup_{x\in\mathbb H^n}
|D^\alpha f_\varepsilon(x)-D^\alpha f(x)|\longrightarrow0.
\]

In either type of chart, the support of the approximant lies in the neighborhood of radius $3\varepsilon$ of the original support, taken in the corresponding coordinate model. For each support, we can choose a slightly larger compact set still contained in the chart domain. Taking $\varepsilon$ sufficiently small ensures that all approximants are supported in these fixed compact sets. They are then transferred to the bundle and extended by zero outside their charts; the resulting sections are smooth and supported in $\operatorname{Int}L$. Equivalence of coordinate and covariant norms on these compact sets, given by Theorem~\ref{norma ck independiente de metrica y conexion}, turns convergence of the components into $C^m$ convergence of the sections. Since there are only finitely many charts and components, choose a common decreasing sequence $\varepsilon_j\to0^{+}$ and sum the approximants. This gives $\mathbf{u}_j\in\Gamma_L(\mathbf{E})$ with $\mathbf{u}_j\to\mathbf{u}$ in $C^m(L)$. \end{proof}

\begin{remark}\label{rem: extension continua distribuciones orden m} Let $\mathbf{F}\colon\mathcal D(M,\mathbf{E}^*)\to W$ be a distribution of global order $m$. For $\mathbf{u}\in\Gamma_K^m(\mathbf{E}^*)$, choose a compact set $L$ with $K\subseteq\operatorname{Int}L$ and, by the preceding lemma, $\mathbf{u}_j\in\Gamma_L(\mathbf{E}^*)$ with $\mathbf{u}_j\to \mathbf{u}$ in $C^m(L)$. The order estimate on $L$ shows that $(\mathbf{F}(\mathbf{u}_j))$ is Cauchy; define \[
 \widetilde{\mathbf{F}}(\mathbf{u}):=\lim_j\mathbf{F}(\mathbf{u}_j).
\] The value is independent of $L$ and the approximation. Indeed, for two choices $L_1,L_2$, take a compact set $L_3$ containing both in its interior. The differences of the two approximants converge to zero in $C^m(L_3)$, and the order estimate on $L_3$ forces the limits to agree. This gives a unique continuous linear extension \[
 \widetilde{\mathbf{F}}\colon\Gamma_c^m(\mathbf{E}^*)\longrightarrow W
\] agreeing with $\mathbf{F}$ on $\mathcal D(M,\mathbf{E}^*)$. For each compact set $K$, its restriction to $\Gamma_K^m(\mathbf{E}^*)$ is continuous because, choosing $L$ as above, $\|\mathbf{u}\|_{C^m(L)}=\|\mathbf{u}\|_{C^m(K)}$. \end{remark}

We give some examples of distributions.

\begin{example}\label{ej:distribuciones-en-haces-continuacion-presentamos-dos-ejemplos-fundamentales-distribuc} We present two fundamental examples of distributions. Let $(M,\mathbf{g})$ be a Riemannian manifold with or without boundary, and let $\mathbf{E}\to M$ be a smooth vector bundle. Fix an auxiliary bundle metric $\mathbf{h}_{\mathbf{E}}$ on $\mathbf{E}$, and denote by $\mathbf{h}_{\mathbf{E}}^*$ the dual metric on $\mathbf{E}^*$.

\begin{enumerate} \item The Dirac delta. Let $p\in M$. Define \[
 \delta_{p}\colon \mathcal{D}(M,\mathbf{E}^{*}) \longrightarrow \mathbf{E}^{*}_{p}, \qquad \delta_{p}(\boldsymbol{\phi}) := \boldsymbol{\phi}(p).
 \] Let $K\subseteq M$ be compact and $\boldsymbol{\phi}\in \mathcal{D}_{K}(M,\mathbf{E}^{*})$. If $p\notin K$, then $\delta_p(\boldsymbol{\phi})=0$. If $p\in K$, we have \[
 |\delta_{p}(\boldsymbol{\phi})|_{\mathbf{h}_{\mathbf{E}}^*(p)} = |\boldsymbol{\phi}(p)|_{\mathbf{h}_{\mathbf{E}}^*(p)} \leq \sup_{x\in K}|\boldsymbol{\phi}(x)|_{\mathbf{h}_{\mathbf{E}}^*(x)} \leq \|\boldsymbol{\phi}\|_{C^{0}(K,\mathbf{E}^*)}.
 \] In either case, there exists $C_K=1$ such that \[
 |\delta_p(\boldsymbol{\phi})|_{\mathbf{h}_{\mathbf{E}}^*(p)} \leq C_K \|\boldsymbol{\phi}\|_{C^0(K,\mathbf{E}^*)}.
 \] By part (d) of Theorem~\ref{teo: equivalencias de continuidad en D'(M,E,W)}, $\delta_p$ is continuous. Moreover, since the estimate holds with $m=0$ for every compact set $K$, we conclude that $\delta_p$ is a distribution of order $0$.

\item Locally integrable sections. Let $\mathbf{f}\in L^{1}_{\operatorname{loc}}(M,\mathbf{E})$. Define \[
 \mathbf{T}_{\mathbf{f}}(\boldsymbol{\phi}) := \int_{M}\boldsymbol{\phi}(\mathbf{f})\, d\lambda_{\mathbf{g}}, \qquad \boldsymbol{\phi}\in \mathcal{D}(M,\mathbf{E}^*).
 \] Let $K\subseteq M$ be compact and $\boldsymbol{\phi}\in \mathcal{D}_{K}(M,\mathbf{E}^{*})$. Then \[
 |\mathbf{T}_{\mathbf{f}}(\boldsymbol{\phi})| \leq \int_{K} |\boldsymbol{\phi}(x)(\mathbf{f}(x))|\, d\lambda_{\mathbf{g}}(x)
 \leq \int_{K} |\boldsymbol{\phi}(x)|_{\mathbf{h}_{\mathbf{E}}^*(x)}\, |\mathbf{f}(x)|_{\mathbf{h}_{\mathbf{E}}(x)}\, d\lambda_{\mathbf{g}}(x).
 \] Hence \[
 |\mathbf{T}_{\mathbf{f}}(\boldsymbol{\phi})| \leq \|\boldsymbol{\phi}\|_{C^{0}(K,\mathbf{E}^*)} \int_{K} |\mathbf{f}(x)|_{\mathbf{h}_{\mathbf{E}}(x)}\, d\lambda_{\mathbf{g}}(x).
 \] Since $\mathbf{f}\in L^1_{\operatorname{loc}}(M,\mathbf{E})$, the integral is finite; denoting it by $C_K$, we obtain \[
 |\mathbf{T}_{\mathbf{f}}(\boldsymbol{\phi})| \leq C_K \|\boldsymbol{\phi}\|_{C^{0}(K,\mathbf{E}^*)}.
\] By part (d) of Theorem~\ref{teo: equivalencias de continuidad en D'(M,E,W)}, $\mathbf{T}_{\mathbf{f}}$ is continuous and, since it depends only on the norm $C^0$, is a distribution of order $0$. \end{enumerate} \end{example} Although every locally integrable section can be viewed as a distribution, we may ask when a distribution is induced by such a section. Because of their significance, these distributions have a special name, which we record in the following definition: \begin{definition}\label{def: distribucion regular}\index{regular distribution on a bundle@regular distribution on a bundle} Let $(M,\mathbf{g})$ be a Riemannian manifold with or without boundary, and let $\mathbf{E}\to M$ be a smooth vector bundle. A distribution $\mathbf{T}\colon \mathcal{D}(M,\mathbf{E}^*)\longrightarrow \mathbb{K}$ is called \textit{regular} if there exists $\mathbf{f}\in L^{1}_{\operatorname{loc}}(M,\mathbf{E})$ such that \[
\mathbf{T}(\boldsymbol{\phi}) = \int_{M}\boldsymbol{\phi}(\mathbf{f})\, d\lambda_{\mathbf{g}} \qquad \forall \boldsymbol{\phi}\in \mathcal{D}(M,\mathbf{E}^{*}),
\] where $\boldsymbol{\phi}(\mathbf{f})(x) = \boldsymbol{\phi}(x)(\mathbf{f}(x))$. \end{definition}

\begin{remark}\label{obs:distribuciones-en-haces-mapeo-lineal-define-operador-espacio-distribuciones} The linear map $\mathbf{f} \mapsto \mathbf{T}_{\mathbf{f}}$ defines an operator from $L^{1}_{\operatorname{loc}}(M,\mathbf{E})$ to the space of distributions. We will later see that this operator is injective. Thus it is customary to identify the section $\mathbf{f}$ with the distribution it induces, writing simply $\mathbf{f}(\boldsymbol{\phi})$ instead of $\mathbf{T}_{\mathbf{f}}(\boldsymbol{\phi})$. \end{remark} The appearance of $\mathbf{E}^*$ in the space of test sections is explained by the canonical pairing between a bundle and its dual: \begin{remark}\label{rem: motivacion dual haces} The choice to define a distribution as a linear operator \[
\mathbf{T}\colon \mathcal{D}(M,\mathbf{E}^*)\longrightarrow W,
\] where $W$ is a finite-dimensional normed vector space, is motivated by the existence of a canonical pairing between a vector bundle and its dual.

If $\mathbf{f}\in \Gamma(\mathbf{E})$ is a sufficiently regular section, then to each $\boldsymbol{\phi}\in \mathcal{D}(M,\mathbf{E}^*)$ we can naturally associate a scalar by fiberwise contraction \[
\boldsymbol{\phi}(\mathbf{f})(x)=\boldsymbol{\phi}(x)(\mathbf{f}(x)).
\] This allows us to define a linear operator \[
\mathbf{T}_{\mathbf{f}}\colon \mathcal{D}(M,\mathbf{E}^*)\longrightarrow \mathbb{K},\qquad
\mathbf{T}_{\mathbf{f}}(\boldsymbol{\phi})
=
\displaystyle\int_M \boldsymbol{\phi}(x)(\mathbf{f}(x))\, d\lambda_{\mathbf{g}}(x),
\] generalizing the classical scalar case.

For example, if $\Omega\subseteq \mathbb{R}^n$ is open and $f\in L^1_{\mathrm{loc}}(\Omega)$, then the expression \[
\displaystyle\int_{\Omega} f(x)\,\varphi(x)\,dx,
\qquad \varphi\in C_c^\infty(\Omega),
\] defines a scalar distribution. The preceding formulation is simply the intrinsic version of this construction in the setting of vector bundles.

\medskip

Observe that this construction uses only the canonical pairing \[
\mathbf{E}_x^* \times \mathbf{E}_x \longrightarrow \mathbb{K},
\] which requires no additional structure. In particular, if we worked with sections of $\mathbf{E}$ instead of sections of $\mathbf{E}^*$, there would be no natural way to define an expression of the form \[
\displaystyle\int_M \boldsymbol{\psi}(x)\, \mathbf{f}(x)\, d\lambda_{\mathbf{g}}(x),
\] for $\boldsymbol{\psi},\mathbf{f}\in \Gamma(\mathbf{E})$ without introducing a bundle metric identifying $\mathbf{E}$ with $\mathbf{E}^*$.

For this reason, the formulation in terms of $\mathbf{E}^*$ is intrinsic and independent of additional geometric choices.

\medskip

The extension to the $W$-valued case arises naturally when considering finite families of distributions. If $W\cong \mathbb{K}^r$, a linear operator \[
\mathbf{T}\colon \mathcal{D}(M,\mathbf{E}^*)\longrightarrow W
\] can be written as \[
\mathbf{T}(\boldsymbol{\phi})=\bigl(T_1(\boldsymbol{\phi}),\dots,T_r(\boldsymbol{\phi})\bigr),
\] where each $T_i\colon \mathcal{D}(M,\mathbf{E}^*)\longrightarrow \mathbb{K}$ is a scalar distribution. In this sense, distributions with values in $W$ introduce no new analytic difficulty, but allow finitely many scalar distributions to be treated simultaneously.

\medskip

This formulation will be fundamental later when considering distributions on product manifolds and describing linear operators by distributional kernels. \end{remark} \subsection{Partial differential operators and distributions} We have seen that distributions are defined as linear functionals on the space of test sections $\mathcal{D}(M,\mathbf{E}^*)$. In particular, to extend the action of differential operators to distributions, we must understand how these operators interact with this space.

Let \[
P\colon \Gamma(\mathbf{E})\longrightarrow \Gamma(\mathbf{F})
\] be a differential operator. If $\mathbf{T}$ is a distribution, we would like to define a new distribution $P\mathbf{T}$ representing the action of $P$ on $\mathbf{T}$. By analogy with the classical case, it is natural to expect this action to be determined by an identity of the form \[
(P\mathbf{T})(\boldsymbol{\phi})=\mathbf{T}(\text{expression depending on }P\text{ and }\boldsymbol{\phi}).
\]

Since $\mathbf{T}$ is defined only on $\mathcal{D}(M,\mathbf{E}^*)$, the preceding expression suggests that to each $\boldsymbol{\phi}\in \mathcal{D}(M,\mathbf{F}^*)$ we should associate a section of $\mathcal{D}(M,\mathbf{E}^*)$. That is, we need to construct an operator \[
P'\colon \mathcal{D}(M,\mathbf{F}^*)\longrightarrow \mathcal{D}(M,\mathbf{E}^*),
\] that transfers the action of $P$ to test sections.

Unlike the pointwise case, transferring differential operators requires the integral structure of the manifold in order to integrate by parts. This suggests that the operator $P'$ should be characterized by an integral identity involving the natural evaluation between bundles and their duals.

\begin{definition}[Formal transpose] \label{def: adjunto formal duales} \label{def:transpuesto-formal-haces} Let $(M,\mathbf{g})$ be a Riemannian manifold without boundary, and let $\mathbf{E},\mathbf{F}\to M$ be complex vector bundles. If $P\in\mathbf{PDO}(\mathbf{E},\mathbf{F})$, its \emph{formal transpose} is an operator \[
P'\colon\Gamma(\mathbf{F}^*)\longrightarrow\Gamma(\mathbf{E}^*)
\] satisfying \begin{equation}
\label{eq:identidad-transpuesto-formal}
\int_M\boldsymbol{\psi}(P\mathbf{u})\,d\lambda_{\mathbf{g}}
=\int_M(P'\boldsymbol{\psi})(\mathbf{u})\,d\lambda_{\mathbf{g}}
\end{equation} for every $\mathbf{u}\in\Gamma_c(\mathbf{E})$ and every $\boldsymbol{\psi}\in\Gamma_c(\mathbf{F}^*)$. The pairing in this identity is complex-bilinear. In the real case, the same definition applies, with all conjugations below omitted. \end{definition}

\begin{proposition}[Intrinsic construction of the formal transpose] \label{prop: dual es PDO integral} \label{prop:transpuesto-formal-riesz-conjugado} Let $(M,\mathbf{g})$ be a Riemannian manifold without boundary, and let $\mathbf{E},\mathbf{F}\to M$ be complex vector bundles. Let $P\in\mathbf{PDO}^{(m)}(\mathbf{E},\mathbf{F})$. The formal transpose exists, is unique, and belongs to $\mathbf{PDO}^{(m)}(\mathbf{F}^*,\mathbf{E}^*)$. For any auxiliary Hermitian metrics $\mathbf{h}_{\mathbf{E}},\mathbf{h}_{\mathbf{F}}$, linear in the first variable, we have \begin{equation}
\label{eq:transpuesto-formal-riesz-conjugado}
P'=\mathcal R_{\mathbf{E}}\circ\overline{P_h^*}\circ\mathcal R_{\mathbf{F}}^{-1}.
\end{equation} The domains and codomains of the composition are displayed in the diagram \[
\mathbf{F}^*
\xrightarrow{\ \mathcal R_{\mathbf{F}}^{-1}\ }
\overline{\mathbf{F}}
\xrightarrow{\ \overline{P_h^*}\ }
\overline{\mathbf{E}}
\xrightarrow{\ \mathcal R_{\mathbf{E}}\ }
\mathbf{E}^*.
\] In particular, the right-hand side of \eqref{eq:transpuesto-formal-riesz-conjugado} is independent of the auxiliary bundle metrics. \end{proposition}

\begin{proof} The Riesz isomorphisms $\mathcal R_{\mathbf{F}}^{-1}\colon \mathbf{F}^*\to\overline{\mathbf{F}}$ and $\mathcal R_{\mathbf{E}}\colon\overline{\mathbf{E}}\to \mathbf{E}^*$ are complex-linear and of order zero. By the theory of operator conjugation, $\overline{P_h^*}\colon\Gamma(\overline{\mathbf{F}})\to\Gamma(\overline{\mathbf{E}})$ is complex-linear and has order at most $m$. By composition, the right-hand side of \eqref{eq:transpuesto-formal-riesz-conjugado} is a complex-linear operator of order at most $m$.

Let $\mathbf{u}\in\Gamma_c(\mathbf{E})$ and $\boldsymbol{\psi}\in\Gamma_c(\mathbf{F}^*)$, and write \[
\mathcal R_{\mathbf{F}}^{-1}\boldsymbol{\psi}=\overline{\mathbf{v}},
\qquad \mathbf{v}\in\Gamma_c(\mathbf{F}).
\] By the definition $\mathcal R_{\mathbf{F}}(\overline{\mathbf{v}})(\mathbf{w})=\mathbf{h}_{\mathbf{F}}(\mathbf{w},\mathbf{v})$ and the Hermitian adjoint identity, \begin{align*}
\int_M\boldsymbol{\psi}(P\mathbf{u})\,d\lambda_{\mathbf{g}}
&=\int_M \mathbf{h}_{\mathbf{F}}(P\mathbf{u},\mathbf{v})\,d\lambda_{\mathbf{g}}\\
&=\int_M \mathbf{h}_{\mathbf{E}}(\mathbf{u},P_h^*\mathbf{v})\,d\lambda_{\mathbf{g}}\\
&=\int_M
\mathcal R_{\mathbf{E}}\!\left(\overline{P_h^*\mathbf{v}}\right)(\mathbf{u})\,d\lambda_{\mathbf{g}}\\
&=\int_M
\bigl(\mathcal R_{\mathbf{E}}\circ\overline{P_h^*}
\circ\mathcal R_{\mathbf{F}}^{-1}\bigr)(\boldsymbol{\psi})(\mathbf{u})\,d\lambda_{\mathbf{g}}.
\end{align*} This proves \eqref{eq:identidad-transpuesto-formal}.

For uniqueness, let $S_1,S_2\in\mathbf{PDO}(\mathbf{F}^*,\mathbf{E}^*)$ be two operators satisfying this identity. Fix $\boldsymbol{\psi}\in\Gamma_c(\mathbf{F}^*)$ and set $\boldsymbol{\eta}=(S_1-S_2)\boldsymbol{\psi}\in\Gamma_c(\mathbf{E}^*)$. If $\mathcal R_{\mathbf{E}}^{-1}\boldsymbol{\eta}=\overline{\mathbf{w}}$, then $\mathbf{w}\in\Gamma_c(\mathbf{E})$ and, taking $\mathbf{u}=\mathbf{w}$, we obtain \[
0=\int_M\boldsymbol{\eta}(\mathbf{w})\,d\lambda_{\mathbf{g}}
=\int_M \mathbf{h}_{\mathbf{E}}(\mathbf{w},\mathbf{w})\,d\lambda_{\mathbf{g}}.
\] Therefore $\mathbf{w}=0$ and $S_1\boldsymbol{\psi}=S_2\boldsymbol{\psi}$. For a section $\boldsymbol{\psi}$ without compact support, choose a cutoff equal to one near a given point and use locality of $S_1,S_2$. Thus $S_1=S_2$ on all sections. Since any choice of auxiliary metrics produces an operator satisfying the same integral identity, uniqueness proves independence of these metrics. \end{proof}

\begin{corollary}[Local formula for the formal transpose] \label{cor:formula-local-transpuesto-formal} Let $(M,\mathbf{g})$ be a Riemannian manifold without boundary, let $\mathbf{E},\mathbf{F}\to M$ be complex vector bundles, and let $P\in\mathbf{PDO}^{(m)}(\mathbf{E},\mathbf{F})$. In a chart $\phi\colon U\to\Omega\subseteq\mathbb R^n$, with coordinate functions $(x^1,\dots,x^n)$, let $(\mathbf{e}_1,\dots,\mathbf{e}_r)$ and $(\mathbf{f}_1,\dots,\mathbf{f}_s)$ be local frames of $\mathbf{E}$ and $\mathbf{F}$, and let $(\mathbf{e}^1,\dots,\mathbf{e}^r)$ and $(\mathbf{f}^1,\dots,\mathbf{f}^s)$ be their dual frames. Set \[
D_\phi^\alpha a
:=
\bigl[D^\alpha(a\circ\phi^{-1})\bigr]\circ\phi.
\] If \[
P\left(\sum_{i=1}^r u^i\mathbf e_i\right)
=
\sum_{j=1}^s
\left(
\sum_{i=1}^r\sum_{|\alpha|\leq m}
a_{\alpha i}^{\,j}D_\phi^\alpha u^i
\right)\mathbf f_j,
\] then, for $\displaystyle \boldsymbol{\psi}=\displaystyle\sum_{j=1}^s\psi_j\mathbf{f}^j$, \begin{equation}
\label{eq:formula-local-transpuesto-formal}
P'\boldsymbol{\psi}
=
\frac{1}{\sqrt{\det(\mathbf{g})}}
\sum_{i=1}^r
\left[
\sum_{j=1}^s\sum_{|\alpha|\leq m}(-1)^{|\alpha|}
D_\phi^\alpha\!\left(
\sqrt{\det(\mathbf{g})}\,
a_{\alpha i}^{\,j}\psi_j
\right)
\right]\mathbf e^i,
\end{equation} where $g_{ab}=\mathbf g(\boldsymbol{\partial}_a,
\boldsymbol{\partial}_b)$. \end{corollary}

\begin{proof} Let $\displaystyle \mathbf u=\displaystyle\sum_{i=1}^{r} u^i\mathbf e_i\in\Gamma_c(\mathbf E|_U)$ and let $\displaystyle \boldsymbol\psi=\displaystyle\sum_{j=1}^{s}\psi_j\mathbf f^j\in\Gamma(\mathbf F^*|_U)$. By the local formulas for $P$ and the Riemannian measure, \begin{align*}
\int_U\boldsymbol\psi(P\mathbf u)\,d\lambda_{\mathbf g}
&=\sum_{i=1}^r\sum_{j=1}^s\sum_{|\alpha|\leq m}
\int_\Omega
\bigl(a_{\alpha i}^{\,j}\circ\phi^{-1}\bigr)
\bigl(\psi_j\circ\phi^{-1}\bigr)
D^\alpha\bigl(u^i\circ\phi^{-1}\bigr)
\bigl(\sqrt{\det(\mathbf{g})}\circ\phi^{-1}\bigr)
\,d\lambda_n.
\end{align*} The function $u^i\circ\phi^{-1}$ has compact support in $\Omega$. Integrating by parts $|\alpha|$ times, without boundary terms, the summand corresponding to $i,j,\alpha$ becomes \[
(-1)^{|\alpha|}
\int_\Omega
\bigl(u^i\circ\phi^{-1}\bigr)
D^\alpha\!\left[
\left(
\sqrt{\det(\mathbf{g})}\,
a_{\alpha i}^{\,j}\psi_j
\right)\circ\phi^{-1}
\right]d\lambda_n.
\] Summing all terms and inserting $\bigl(\sqrt{\det(\mathbf{g})}\bigr)^{-1}
\sqrt{\det(\mathbf{g})}$ gives \[
\int_U\boldsymbol\psi(P\mathbf u)\,d\lambda_{\mathbf g}
=
\int_U(P'\boldsymbol\psi)(\mathbf u)\,d\lambda_{\mathbf g},
\] where $P'\boldsymbol\psi$ is precisely the right-hand side of \eqref{eq:formula-local-transpuesto-formal}. Uniqueness of the formal transpose identifies this local expression with the global operator already constructed. \end{proof}

\begin{proposition}[Rules for the formal transpose] \label{prop:reglas-transpuesto-formal} Let $(M,\mathbf{g})$ be a Riemannian manifold without boundary, and let $\mathbf{E},\mathbf{F},\mathbf{G}\to M$ be complex vector bundles. Let $P\in\mathbf{PDO}(\mathbf{E},\mathbf{F})$, $Q\in\mathbf{PDO}(\mathbf{F},\mathbf{G})$, and $R\in\mathbf{PDO}(\mathbf{E},\mathbf{F})$. Then: \begin{enumerate}[label=(\alph*)] \item for every $\lambda\in\mathbb C$, \[
(P+\lambda R)'=P'+\lambda R';
\] \item \[
(Q\circ P)'=P'\circ Q';
\] \item if $U\subseteq M$ is open and $P_U$ denotes the local restriction of $P$, then \[
(P_U)'=(P')_U;
\] \item if $P$ has order zero, then \[
(P'\boldsymbol{\psi})(\mathbf{u})=\boldsymbol{\psi}(P\mathbf{u})
\] pointwise; in particular, $(M_f)'=M_f$; \item if \[
\kappa_{\mathbf{E}}\colon\overline{\mathbf{E}^*}\longrightarrow(\overline{\mathbf{E}})^*,
\qquad
\kappa_{\mathbf{E}}(\overline{\boldsymbol{\eta}})(\overline v)
:=\overline{\boldsymbol{\eta}(v)},
\] and \[
\kappa_{\mathbf{F}}\colon\overline{\mathbf{F}^*}\longrightarrow(\overline{\mathbf{F}})^*,
\qquad
\kappa_{\mathbf{F}}(\overline\psi)(\overline w)
:=\overline{\psi(w)},
\] then \begin{equation}
\label{eq:transpuesto-conjugacion}
(\overline P)'
=\kappa_{\mathbf{E}}\circ\overline{P'}\circ\kappa_{\mathbf{F}}^{-1}.
\end{equation} \end{enumerate} \end{proposition}

\begin{proof} The first identity follows from bilinearity of \eqref{eq:identidad-transpuesto-formal}; unlike the Hermitian adjoint, no $\overline\lambda$ appears here. For $\mathbf{u}\in\Gamma_c(\mathbf{E})$ and $\boldsymbol{\theta}\in\Gamma_c(\mathbf{G}^*)$, \[
\int_M\boldsymbol{\theta}(QP\mathbf{u})\,d\lambda_{\mathbf{g}}
=\int_M(Q'\boldsymbol{\theta})(P\mathbf{u})\,d\lambda_{\mathbf{g}}
=\int_M(P'Q'\boldsymbol{\theta})(\mathbf{u})\,d\lambda_{\mathbf{g}},
\] which proves the second identity by uniqueness. The same integral identity, applied entirely within $U$, proves the third. The fourth is the definition of the fiberwise algebraic transpose.

The map $\kappa_{\mathbf{E}}$ is a complex-linear isomorphism. For $\mathbf{u}\in\Gamma_c(\mathbf{E})$ and $\boldsymbol{\eta}\in\Gamma_c(\mathbf{F}^*)$, we have \[
\kappa_{\mathbf{F}}(\overline{\boldsymbol{\eta}})(\overline P\,\overline{\mathbf{u}})
=\overline{\boldsymbol{\eta}(P\mathbf{u})},
\qquad
\kappa_{\mathbf{E}}(\overline{P'\boldsymbol{\eta}})(\overline{\mathbf{u}})
=\overline{(P'\boldsymbol{\eta})(\mathbf{u})}.
\] Conjugating the integral identity for $P'$ yields the integral identity for the right-hand side of \eqref{eq:transpuesto-conjugacion}; uniqueness of the formal transpose finishes the proof. \end{proof} The preceding construction extends the action of differential operators to distributions. If $P\colon \Gamma(\mathbf{E})\longrightarrow \Gamma(\mathbf{F})$ is a differential operator and $P'\colon \Gamma(\mathbf{F}^*)\longrightarrow \Gamma(\mathbf{E}^*)$ denotes its formal transpose with respect to the bilinear evaluation pairing, then $P'$ sends test sections of $\mathbf{F}^*$ to test sections of $\mathbf{E}^*$. This allows the action of $P$ on a distribution to be defined by duality.

\begin{lemma}[$C^k$ estimate for differential operators]\label{lem: estimacion Ck pdo compacto}\index{estimate for differential operators@estimate for differential operators} Let $\mathbf{E},\mathbf{F}\longrightarrow M$ be smooth vector bundles of ranks $r$ and $s$, respectively, over a smooth manifold $M$ with or without boundary, and let $P\in \mathbf{PDO}^{(m)}(\mathbf{E},\mathbf{F})$. Then, for every compact set $K\subseteq M$ and every $k\in\mathbb{N}_{0}$, there exists a constant $C_{K,k}>0$ such that \[
\|P\mathbf{u}\|_{C^k(K,\mathbf{F})}
\leq
C_{K,k}\|\mathbf{u}\|_{C^{k+m}(K,\mathbf{E})},
\qquad
\mathbf{u}\in\Gamma^{k+m}_{K}(\mathbf{E}).
\] In particular, $P$ induces a continuous map \[
P\colon \mathcal{D}_{K}(M,\mathbf{E})\longrightarrow \mathcal{D}_{K}(M,\mathbf{F}).
\] \end{lemma}

\begin{proof} Let $K\subseteq M$ be compact. Take a finite family of smooth charts $(U_i,\phi_i)$, with $i\in\{1,\dots,N\}$, such that each $U_i$ is a regular coordinate ball, local trivializations \[
\boldsymbol{\tau}_{\mathbf{E},i}\colon \mathbf{E}\restriction_{U_i}\longrightarrow U_i\times \mathbb{K}^{r},
\qquad
\boldsymbol{\tau}_{\mathbf{F},i}\colon \mathbf{F}\restriction_{U_i}\longrightarrow U_i\times \mathbb{K}^{s},
\] exist, and $K\subseteq \displaystyle\bigcup_{i=1}^{N}U_i$ holds.

Consider a partition of unity $\{\psi_i\}_{i=1}^{N}$ subordinate to the cover $\{U_i\}_{i=1}^{N}$ of $K$, with sum one on a neighborhood of $K$, and define \[
K_i:=K\cap\operatorname{supp}(\psi_i).
\] Since $\operatorname{supp}(\psi_i)$ is closed in $M$ and contained in $U_i$, each $K_i$ is a compact subset of $U_i$. Omit the indices for which $K_i=\varnothing$ and renumber the rest. If $\mathbf{u}\in\Gamma^{k+m}_{K}(\mathbf{E})$, then $\supp(\mathbf{u})\subseteq K$ and \[
\mathbf{u}=\displaystyle\sum_{i=1}^{N}\psi_i \mathbf{u}.
\] Moreover, since $\supp(\psi_i)\subseteq U_i$, we have \[
\supp(\psi_i \mathbf{u})\subseteq \supp(\mathbf{u})\cap \supp(\psi_i)
\subseteq K\cap\operatorname{supp}(\psi_i)
=
K_i.
\]

Fix $i\in\{1,\dots,N\}$ and write $\mathbf{v}_i:=\psi_i \mathbf{u}$. Define its local representation in the trivialization $\boldsymbol{\tau}_{\mathbf{E},i}$ by \[
\widetilde{v}_i:=\pi_2\circ \boldsymbol{\tau}_{\mathbf{E},i}\circ (\mathbf{v}_i\restriction_{U_i})\colon U_i\longrightarrow \mathbb{K}^{r},
\] and similarly define the local representation of $P\mathbf{v}_i$ in the trivialization $\boldsymbol{\tau}_{\mathbf{F},i}$ by \[
\widetilde{P\mathbf{v}_i}:=\pi_2\circ \boldsymbol{\tau}_{\mathbf{F},i}\circ ((P\mathbf{v}_i)\restriction_{U_i})\colon U_i\longrightarrow \mathbb{K}^{s}.
\]

By Theorem~\ref{thm:pdo expresion local explicita}, there exist smooth functions \[
a_{\alpha,i}\colon U_i\longrightarrow \operatorname{Hom}(\mathbb{K}^{r},\mathbb{K}^{s}),
\qquad |\alpha|\leq m,
\] such that \[
\widetilde{P\mathbf{v}_i}
=
\displaystyle\sum_{|\alpha|\leq m}a_{\alpha,i}\,D^\alpha \widetilde{v}_i.
\]

Let $\beta$ be a multi-index with $|\beta|\leq k$. We apply the Leibniz rule to the bilinear map \[
\operatorname{Hom}(\mathbb{K}^{r},\mathbb{K}^{s})\times \mathbb{K}^{r}
\longrightarrow
\mathbb{K}^{s},
\qquad
(A,w)\longmapsto A(w).
\] To make the notation precise, if $A\colon U_i\longrightarrow \operatorname{Hom}(\mathbb{K}^{r},\mathbb{K}^{s})$ and we write $A(x)=(A^{b}_{a}(x))$, with $1\leq a\leq r$ and $1\leq b\leq s$, then we define $D^\gamma A$ componentwise, that is, \[
D^\gamma A(x):=\bigl(D^\gamma A^{b}_{a}(x)\bigr).
\] With this convention, $D^\gamma A(x)$ remains an element of $\operatorname{Hom}(\mathbb{K}^{r},\mathbb{K}^{s})$.

Then, for each $|\beta|\leq k$, we have \[
D^\beta(\widetilde{P\mathbf{v}_i})
=
\displaystyle\sum_{|\alpha|\leq m}
\displaystyle\sum_{\gamma\leq \beta}
{\beta\choose\gamma}
\bigl(D^\gamma a_{\alpha,i}\bigr)
\bigl(D^{\alpha+\beta-\gamma}\widetilde{v}_i\bigr).
\] This identity follows by applying the Leibniz rule componentwise to the preceding matrix expression.

Since $K_i$ is compact and the coefficients $a_{\alpha,i}$ are smooth, the functions $D^\gamma a_{\alpha,i}$ are bounded on $K_i$ for every $|\gamma|\leq k$. Moreover, if $|\alpha|\leq m$, $\gamma\leq\beta$, and $|\beta|\leq k$, then $|\alpha+\beta-\gamma|\leq m+k$. Thus there exists a constant $A_i>0$ such that \[
\max_{|\beta|\leq k}\|D^\beta(\widetilde{P\mathbf v_i})\|_{C^0(K_i,\mathbb K^s)}
\leq
A_i
\max_{|\eta|\leq k+m}\|D^\eta\widetilde v_i\|_{C^0(K_i,\mathbb K^r)}.
\]

Now relate this estimate to the $C^j$ norms. On $U_i$, write $\phi_i=(x_i^1,\dots,x_i^n)$ and denote the coordinate vector fields by $\boldsymbol{\partial}_{x_i^1},\dots,\boldsymbol{\partial}_{x_i^n}$. Define a connection $\widetilde{\nabla}^{TM,i}$ on $TU_i$ by declaring these coordinate vector fields parallel, that is, \[
\widetilde{\nabla}^{TM,i}_{\boldsymbol{\partial}_{x_i^\ell}}\boldsymbol{\partial}_{x_i^a}=0,
\qquad
a,\ell\in\{1,\dots,n\}.
\] Likewise, if $\mathbf{e}_1^i,\dots,\mathbf{e}_r^i$ is the local frame of $\mathbf{E}\restriction_{U_i}$ induced by the trivialization $\boldsymbol{\tau}_{\mathbf{E},i}$, define a connection $\widetilde{\nabla}^{\mathbf{E},i}$ on $\mathbf{E}\restriction_{U_i}$ by declaring this frame parallel: \[
\widetilde{\nabla}^{\mathbf{E},i}_{\boldsymbol{\partial}_{x_i^\ell}}\mathbf{e}_a^i=0,
\qquad
a\in\{1,\dots,r\},\ \ell\in\{1,\dots,n\}.
\] Also take the Riemannian metric $\widetilde{\mathbf{g}}_i$ making the coordinate frame orthonormal and the bundle metric $\widetilde{\mathbf{h}}_{\mathbf{E},i}$ making the frame $\mathbf{e}_1^i,\dots,\mathbf{e}_r^i$ orthonormal.

With these choices, if $\mathbf{v}_i=\displaystyle\sum_{a=1}^{r}v_i^a\mathbf{e}_a^i$ on $U_i$, then $\widetilde{v}_i=(v_i^1,\dots,v_i^r)$ and the components of the iterated covariant derivatives of $\mathbf{v}_i$, computed using the connections $\widetilde{\nabla}^{TM,i}$ and $\widetilde{\nabla}^{\mathbf{E},i}$, are precisely the partial derivatives of the components of $\widetilde{v}_i$. Consequently, for each $j\leq k+m$ there exists a constant $B_{i,j}>0$ such that \[
\max_{|\eta|\leq j}\|D^\eta\widetilde v_i\|_{C^0(K_i,\mathbb K^r)}
\leq
B_{i,j}\|\mathbf{v}_i\|_{\widetilde{C^j}(K_i,\mathbf{E})}.
\] Here $\|\cdot\|_{\widetilde{C^j}(K_i,\mathbf{E})}$ denotes the $C^j$ norm on sections of $\mathbf{E}$ constructed using the choices $\widetilde{\mathbf{g}}_i$, $\widetilde{\nabla}^{TM,i}$, $\widetilde{\mathbf{h}}_{\mathbf{E},i}$, and $\widetilde{\nabla}^{\mathbf{E},i}$.

By Theorem~\ref{norma ck independiente de metrica y conexion}, applied to the compact set $K_i$ and the bundle $\mathbf{E}$, the norms $\|\cdot\|_{\widetilde{C^{k+m}}(K_i,\mathbf{E})}$ and $\|\cdot\|_{C^{k+m}(K_i,\mathbf{E})}$ are equivalent. In particular, there exists a constant $C_i>0$ such that \[
\|\mathbf{v}_i\|_{\widetilde{C^{k+m}}(K_i,\mathbf{E})}
\leq
C_i\|\mathbf{v}_i\|_{C^{k+m}(K_i,\mathbf{E})}.
\]

Since $\mathbf{v}_i=\psi_i \mathbf{u}$, the Leibniz rule for covariant derivatives implies the existence of a constant $L_i>0$ such that \[
\|\mathbf{v}_i\|_{C^{k+m}(K_i,\mathbf{E})}
=
\|\psi_i \mathbf{u}\|_{C^{k+m}(K_i,\mathbf{E})}
\leq
L_i\|\mathbf{u}\|_{C^{k+m}(K,\mathbf{E})}.
\]

Apply the same reasoning to the bundle $\mathbf{F}$. If $\mathbf{f}_1^i,\dots,\mathbf{f}_s^i$ is the local frame of $\mathbf{F}\restriction_{U_i}$ induced by $\boldsymbol{\tau}_{\mathbf{F},i}$, take the connection $\widetilde{\nabla}^{\mathbf{F},i}$ determined by \[
\widetilde{\nabla}^{\mathbf{F},i}_{\boldsymbol{\partial}_{x_i^\ell}}\mathbf{f}_b^i=0,
\qquad
b\in\{1,\dots,s\},\ \ell\in\{1,\dots,n\},
\] and the bundle metric $\widetilde{h}_{\mathbf{F},i}$ making that frame orthonormal. Then, using Theorem~\ref{norma ck independiente de metrica y conexion} again, there exists a constant $D_i>0$ such that \[
\|P\mathbf{v}_i\|_{C^k(K_i,\mathbf{F})}
\leq
D_i
\max_{|\beta|\leq k}\|D^\beta(\widetilde{P\mathbf v_i})\|_{C^0(K_i,\mathbb K^s)}.
\]

Combining the preceding estimates gives \[
\|P\mathbf{v}_i\|_{C^k(K_i,\mathbf{F})}
\leq
D_iA_iB_{i,k+m}C_iL_i\|\mathbf{u}\|_{C^{k+m}(K,\mathbf{E})}.
\] Thus, for each $i\in\{1,\dots,N\}$ there exists a constant $C_i'>0$ such that \[
\|P(\psi_i \mathbf{u})\|_{C^k(K_i,\mathbf{F})}
\leq
C_i'\|\mathbf{u}\|_{C^{k+m}(K,\mathbf{E})}.
\]

Since $\mathbf{u}=\displaystyle\sum_{i=1}^{N}\psi_i \mathbf{u}$ and $P$ is linear, we have $P\mathbf{u}=\displaystyle\sum_{i=1}^{N}P(\psi_i \mathbf{u})$. Moreover, locality of $P$ implies $\supp(P(\psi_i \mathbf{u}))\subseteq \supp(\psi_i \mathbf{u})\subseteq K_i$. Therefore $\|P(\psi_i \mathbf{u})\|_{C^k(K,\mathbf{F})}=\|P(\psi_i \mathbf{u})\|_{C^k(K_i,\mathbf{F})}$. Using the triangle inequality for the norm $C^k(K,\mathbf{F})$, we obtain \[
\|P\mathbf{u}\|_{C^k(K,\mathbf{F})}
\leq
\displaystyle\sum_{i=1}^{N}\|P(\psi_i \mathbf{u})\|_{C^k(K,\mathbf{F})}
=
\displaystyle\sum_{i=1}^{N}\|P(\psi_i \mathbf{u})\|_{C^k(K_i,\mathbf{F})}.
\] The preceding estimates give \[
\|P\mathbf{u}\|_{C^k(K,\mathbf{F})}
\leq
\left(\displaystyle\sum_{i=1}^{N}C_i'\right)\|\mathbf{u}\|_{C^{k+m}(K,\mathbf{E})}.
\] Setting $C_{K,k}:=\displaystyle\sum_{i=1}^{N}C_i'$ yields the desired estimate.

Locality of $P$ ensures that $\supp(P\mathbf{u})\subseteq \supp(\mathbf{u})$, so the map is well defined, and the preceding argument also shows that it is continuous. \end{proof}

\begin{proposition}\label{prop: pdo induce operador distribuciones} Let $(M,\mathbf{g})$ be a Riemannian manifold without boundary, let $\mathbf{E},\mathbf{F}\longrightarrow M$ be smooth vector bundles, let $W$ be a finite-dimensional $\mathbb{K}$-vector space, and let $P\in \mathbf{PDO}^{(m)}(\mathbf{E},\mathbf{F})$. Then $P$ induces a continuous linear operator \[
P\colon \mathcal{D}'(M,\mathbf{E},W)\longrightarrow \mathcal{D}'(M,\mathbf{F},W)
\] given by \[
(P\mathbf{T})(\boldsymbol{\phi}):=\mathbf{T}(P'\boldsymbol{\phi}),
\qquad
\boldsymbol{\phi}\in \mathcal{D}(M,\mathbf{F}^*),
\] where $P'\colon \Gamma(\mathbf{F}^*)\longrightarrow \Gamma(\mathbf{E}^*)$ is the bilinear formal transpose. Continuity is understood with respect to the weak-$*$ topologies on both distribution spaces. \end{proposition}

\begin{proof} First check that the definition makes sense. By Proposition~\ref{prop: dual es PDO integral}, we have $P'\in \mathbf{PDO}^{(m)}(\mathbf{F}^*,\mathbf{E}^*)$. Thus, if $\boldsymbol{\phi}\in\mathcal{D}(M,\mathbf{F}^*)$, then $P'\boldsymbol{\phi}$ is a smooth section of $\mathbf{E}^*$.

Moreover, since differential operators are local, we have \[
\supp(P'\boldsymbol{\phi})\subseteq \supp(\boldsymbol{\phi}).
\] In particular, if $\boldsymbol{\phi}$ has compact support, then $P'\boldsymbol{\phi}$ also has compact support. Thus $P'\boldsymbol{\phi}\in\mathcal{D}(M,\mathbf{E}^*)$, and consequently the expression $\mathbf{T}(P'\boldsymbol{\phi})$ is well defined.

Define $(P\mathbf{T})(\boldsymbol{\phi}):=\mathbf{T}(P'\boldsymbol{\phi})$ for $\boldsymbol{\phi}\in\mathcal{D}(M,\mathbf{F}^*)$. Linearity of $P\mathbf{T}$ follows from linearity of $\mathbf{T}$ and $P'$.

Now prove that $P\mathbf{T}\in\mathcal{D}'(M,\mathbf{F},W)$. Let $K\subseteq M$ be compact. If $\boldsymbol{\phi}\in\mathcal{D}_{K}(M,\mathbf{F}^*)$, then locality of $P'$ gives $P'\boldsymbol{\phi}\in\mathcal{D}_{K}(M,\mathbf{E}^*)$. Since $\mathbf{T}\in\mathcal{D}'(M,\mathbf{E},W)$, Theorem~\ref{teo: equivalencias de continuidad en D'(M,E,W)} provides $\ell\in\mathbb{N}_{0}$ and $c_K>0$ such that \[
|\mathbf{T}(\boldsymbol{\psi})|_W\leq c_K\|\boldsymbol{\psi}\|_{C^\ell(K,\mathbf{E}^*)}
\qquad
\forall \boldsymbol{\psi}\in\mathcal{D}_{K}(M,\mathbf{E}^*).
\] Applying this estimate to $\boldsymbol{\psi}=P'\boldsymbol{\phi}$ gives \[
|(P\mathbf{T})(\boldsymbol{\phi})|_W
=
|\mathbf{T}(P'\boldsymbol{\phi})|_W
\leq
c_K\|P'\boldsymbol{\phi}\|_{C^\ell(K,\mathbf{E}^*)}.
\] Since $P'$ is a differential operator of order $m$, Lemma~\ref{lem: estimacion Ck pdo compacto} provides a constant $C_{K,\ell}>0$ such that \[
\|P'\boldsymbol{\phi}\|_{C^\ell(K,\mathbf{E}^*)}
\leq
C_{K,\ell}\|\boldsymbol{\phi}\|_{C^{\ell+m}(K,\mathbf{F}^*)}.
\] Therefore \[
|(P\mathbf{T})(\boldsymbol{\phi})|_W
\leq
c_K C_{K,\ell}\|\boldsymbol{\phi}\|_{C^{\ell+m}(K,\mathbf{F}^*)}
\qquad
\forall \boldsymbol{\phi}\in\mathcal{D}_{K}(M,\mathbf{F}^*).
\] Applying Theorem~\ref{teo: equivalencias de continuidad en D'(M,E,W)} again, now to the bundle $\mathbf{F}$, we conclude that $P\mathbf{T}\in\mathcal{D}'(M,\mathbf{F},W)$.

Next check linearity of the operator induced by $P$. If $\lambda,\mu\in\mathbb{K}$ and $\mathbf{T}_1,\mathbf{T}_2\in\mathcal{D}'(M,\mathbf{E},W)$, then, for every $\boldsymbol{\phi}\in\mathcal{D}(M,\mathbf{F}^*)$, \[
P(\lambda \mathbf{T}_1+\mu \mathbf{T}_2)(\boldsymbol{\phi})
=
(\lambda \mathbf{T}_1+\mu \mathbf{T}_2)(P'\boldsymbol{\phi})
=
\lambda \mathbf{T}_1(P'\boldsymbol{\phi})+\mu \mathbf{T}_2(P'\boldsymbol{\phi}).
\] Therefore $P(\lambda \mathbf{T}_1+\mu \mathbf{T}_2)=\lambda P\mathbf{T}_1+\mu P\mathbf{T}_2$.

Finally, prove continuity of the induced operator in the weak-$*$ topologies. Let $(\mathbf{T}_\alpha)_\alpha$ be a net in $\mathcal{D}'(M,\mathbf{E},W)$ such that $\mathbf{T}_\alpha\to \mathbf{T}$ in $\mathcal{D}'(M,\mathbf{E},W)$. By Remark~\ref{obs: top debil * distribuciones haces}, this means that \[
\mathbf{T}_\alpha(\boldsymbol{\psi})\longrightarrow \mathbf{T}(\boldsymbol{\psi})
\qquad
\forall \boldsymbol{\psi}\in\mathcal{D}(M,\mathbf{E}^*).
\] Now let $\boldsymbol{\phi}\in\mathcal{D}(M,\mathbf{F}^*)$. As already shown, $P'\boldsymbol{\phi}\in\mathcal{D}(M,\mathbf{E}^*)$. Therefore \[
(P\mathbf{T}_\alpha)(\boldsymbol{\phi})
=
\mathbf{T}_\alpha(P'\boldsymbol{\phi})
\longrightarrow
\mathbf{T}(P'\boldsymbol{\phi})
=
(P\mathbf{T})(\boldsymbol{\phi}).
\] Since this holds for every $\boldsymbol{\phi}\in\mathcal{D}(M,\mathbf{F}^*)$, Remark~\ref{obs: top debil * distribuciones haces} again gives $P\mathbf{T}_\alpha\to P\mathbf{T}$ in $\mathcal{D}'(M,\mathbf{F},W)$. Thus the operator induced by $P$ is continuous. \end{proof}

\paragraph{Differential operators up to the boundary.}
If $M$ has boundary, define the transpose first on $\operatorname{Int}M$.
In a boundary chart and local frames, write
$P=\sum_{|\alpha|\leq m}a_\alpha\partial^\alpha$ and
$d\lambda_{\mathbf g}=\rho\,dx$, where $\rho>0$ is smooth up to the boundary.
Its local transpose is
\[
 P'\boldsymbol\phi
 =\rho^{-1}\sum_{|\alpha|\leq m}(-1)^{|\alpha|}
       \partial^\alpha(\rho a_\alpha^{\mathsf T}\boldsymbol\phi).
\]
All coefficients extend smoothly to the boundary. These local formulas agree
on overlaps because they agree on the dense interior, where the formal
transpose is unique. They therefore define a differential operator on $M$.
The support and seminorm estimates in the preceding proof remain valid for
test sections smooth up to the boundary. Consequently
$(P\mathbf T)(\boldsymbol\phi):=\mathbf T(P'\boldsymbol\phi)$ defines the same
continuous distributional action on a manifold with boundary.

Here the transpose is an identity of differential expressions; integration
by parts with tests nonzero at the boundary may produce boundary terms.
Thus, for a smooth section $\mathbf u$, $P\mathbf T_{\mathbf u}$ agrees with
the regular distribution of $P\mathbf u$ on the interior, but may also have
a boundary-supported part. For example, on $[0,\infty)$,
$(\partial_x T_1)(\phi)=-\int_0^\infty\phi'(x)\,dx=\phi(0)$.
This convention includes distributions supported on the boundary and does
not impose a boundary condition on $\mathbf u$.

\begin{remark}\label{obs:distribuciones-en-haces-seccion-suave-distribucion-regular-asociada-coincide} Assume again that $M$ has no boundary, as in the preceding proposition. If $\mathbf{u}\in\Gamma(\mathbf{E})$ is a smooth section, then $P\mathbf{u}\in\Gamma(\mathbf{F})$, and the regular distribution associated with $P\mathbf{u}$ agrees with the distribution obtained by applying $P$ to the regular distribution associated with $\mathbf{u}$. If $\mathbf{T}_{\mathbf{u}}$ is given by \[
\mathbf{T}_{\mathbf{u}}(\boldsymbol{\psi}):=\int_M \boldsymbol{\psi}(\mathbf{u})\,d\lambda_{\mathbf{g}},\qquad \boldsymbol{\psi}\in\mathcal{D}(M,\mathbf{E}^*),
\] then, for every $\boldsymbol{\phi}\in\mathcal{D}(M,\mathbf{F}^*)$, \[
(P\mathbf{T}_{\mathbf{u}})(\boldsymbol{\phi})=\mathbf{T}_{\mathbf{u}}(P'\boldsymbol{\phi})=\int_M (P'\boldsymbol{\phi})(\mathbf{u})\,d\lambda_{\mathbf{g}}.
\] By the definition of the bilinear formal transpose, \[
\int_M (P'\boldsymbol{\phi})(\mathbf{u})\,d\lambda_{\mathbf{g}}=\int_M \boldsymbol{\phi}(P\mathbf{u})\,d\lambda_{\mathbf{g}}.
\] Therefore $P\mathbf{T}_{\mathbf{u}}=\mathbf{T}_{P\mathbf{u}}$.

For operators of order $0$, the identity also holds when $M$ has boundary. In that case, $P$ corresponds to a smooth section of $\operatorname{Hom}(\mathbf{E},\mathbf{F})$ and $P'$ corresponds to the pointwise transpose $\mathbf{F}^*\to \mathbf{E}^*$. In the special case where $\mathbf{E}=\mathbf{F}$ and $P$ is multiplication by a function $f\in C^\infty(M)$, we obtain \[
(fT)(\boldsymbol{\phi})=T(f\boldsymbol{\phi}).
\] \end{remark}

The preceding proposition shows that the action of a differential operator on distributions is determined by the action of its bilinear formal transpose on test sections. This construction suffices to define distributional derivatives and, more generally, to transfer differential operators to the distributional setting. We now give some examples illustrating what has been developed so far. \begin{example}[Differential operators and their distributional extensions on bundles]\label{ej:distribuciones-en-haces-operadores-diferenciales-y-sus-extensiones-distribucionales-e} Let $(M,\mathbf{g})$ be a closed Riemannian manifold, let $\mathbf{E},\mathbf{F}\longrightarrow M$ be smooth vector bundles, and let $W$ be a finite-dimensional $\mathbb{K}$-vector space.

If $P\in \mathbf{PDO}^{(m)}(\mathbf{E},\mathbf{F})$, the induced operator on distributions is defined by \[
(P\mathbf{T})(\boldsymbol{\phi}):=\mathbf{T}(P'\boldsymbol{\phi}),
\qquad
\mathbf{T}\in\mathcal{D}'(M,\mathbf{E},W),\ \boldsymbol{\phi}\in\mathcal{D}(M,\mathbf{F}^*),
\] where $P'\colon \Gamma(\mathbf{F}^*)\longrightarrow \Gamma(\mathbf{E}^*)$ is the bilinear formal transpose. This definition is valid for arbitrary distributions with values in $W$.

To make the integral formulas concrete, we also consider the special case $W=\mathbb{K}$. If $\mathbf{u}\in\Gamma(\mathbf{E})$, denote by $\mathbf{T}_{\mathbf{u}}$ the regular distribution induced by $\mathbf{u}$, given by \[
\mathbf{T}_{\mathbf{u}}(\boldsymbol{\psi}):=\int_M \boldsymbol{\psi}(\mathbf{u})\,d\lambda_{\mathbf{g}},
\qquad
\boldsymbol{\psi}\in\mathcal{D}(M,\mathbf{E}^*).
\]

\begin{enumerate}

\item Let $\nabla^{\mathbf{E}}\colon \Gamma(\mathbf{E})\longrightarrow \Gamma(T^*M\otimes \mathbf{E})$ be a connection. The induced operator is \[
\nabla^{\mathbf{E}}\colon \mathcal{D}'(M,\mathbf{E},W)\longrightarrow \mathcal{D}'(M,T^*M\otimes \mathbf{E},W),
\] and is given by \[
(\nabla^{\mathbf{E}}\mathbf{T})(\boldsymbol{\phi})
=
\mathbf{T}\bigl((\nabla^{\mathbf{E}})'\boldsymbol{\phi}\bigr),
\qquad
\boldsymbol{\phi}\in\mathcal{D}(M,TM\otimes \mathbf{E}^*).
\] In the regular case, that is, when $W=\mathbb{K}$ and $\mathbf{T}=\mathbf{T}_{\mathbf{u}}$, we have \[
(\nabla^{\mathbf{E}}\mathbf{T}_{\mathbf{u}})(\boldsymbol{\phi})
=
\int_M \boldsymbol{\phi}(\nabla^{\mathbf{E}}\mathbf{u})\,d\lambda_{\mathbf{g}}.
\]

If $\mathbf{E}$ is also equipped with a real bundle metric in the real case or a Hermitian metric in the complex case, the Riesz isomorphism of the bundle $\mathbf{H}:=T^*M\otimes \mathbf{E}$ uniquely determines $\boldsymbol{\phi}^{\sharp}\in\Gamma_c(\mathbf{H})$ by \[
 \overline{\boldsymbol{\phi}^{\sharp}}=\mathcal R_{\mathbf{H}}^{-1}\boldsymbol{\phi},
 \qquad \boldsymbol{\phi}(w)=\mathbf{h}_{\mathbf{H}}(w,\boldsymbol{\phi}^{\sharp}).
\] In the complex case, $\boldsymbol{\phi}\mapsto\boldsymbol{\phi}^{\sharp}$ is antilinear; the map $\mathcal R_{\mathbf{H}}^{-1}\colon \mathbf{H}^*\to\overline{\mathbf{H}}$ is complex-linear. Then \[
(\nabla^{\mathbf{E}}\mathbf{T}_{\mathbf{u}})(\boldsymbol{\phi})
=
\int_M \mathbf{h}_{\mathbf{H}}(\nabla^{\mathbf{E}}\mathbf{u},\boldsymbol{\phi}^{\sharp})\,d\lambda_{\mathbf{g}}.
\] By the definition of the formal Hermitian adjoint, \[
\int_M \mathbf{h}_{\mathbf{H}}(\nabla^{\mathbf{E}}\mathbf{u},\boldsymbol{\phi}^{\sharp})\,d\lambda_{\mathbf{g}}
=
\int_M \mathbf{h}_{\mathbf{E}}(\mathbf{u},(\nabla^{\mathbf{E}})_h^*\boldsymbol{\phi}^{\sharp})\,d\lambda_{\mathbf{g}}.
\] Thus, in the presence of bundle metrics, the same identity can be read as an integration-by-parts formula.

\item For the covariant derivative of order $k$, \[
\nabla^k\colon \Gamma(\mathbf{E})\longrightarrow \Gamma(T^{(0,k)}(TM)\otimes \mathbf{E}),
\] the induced operator is given by \[
(\nabla^k \mathbf{T})(\boldsymbol{\phi})
=
\mathbf{T}\bigl((\nabla^k)'\boldsymbol{\phi}\bigr),
\qquad
\boldsymbol{\phi}\in\mathcal{D}(M,T^{(k,0)}(TM)\otimes \mathbf{E}^*).
\] If $W=\mathbb{K}$ and $\mathbf{T}=\mathbf{T}_{\mathbf{u}}$, then \[
(\nabla^k\mathbf{T}_{\mathbf{u}})(\boldsymbol{\phi})
=
\int_M \boldsymbol{\phi}(\nabla^k \mathbf{u})\,d\lambda_{\mathbf{g}}.
\] In the presence of a real or Hermitian bundle metric, set $\mathbf{H}_k:=T^{(0,k)}(TM)\otimes \mathbf{E}$ and define $\boldsymbol{\phi}^{\sharp}$ by $\overline{\boldsymbol{\phi}^{\sharp}}=\mathcal R_{\mathbf{H}_k}^{-1}\boldsymbol{\phi}$. We obtain \[
(\nabla^k\mathbf{T}_{\mathbf{u}})(\boldsymbol{\phi})
=
\int_M \mathbf{h}_{\mathbf{H}_k}(\nabla^k \mathbf{u},\boldsymbol{\phi}^{\sharp})\,d\lambda_{\mathbf{g}}
=
\int_M \mathbf{h}_{\mathbf{E}}(\mathbf{u},(\nabla^k)_h^*\boldsymbol{\phi}^{\sharp})\,d\lambda_{\mathbf{g}}.
\]

\item Let $\operatorname{div}\colon \mathfrak{X}(M)\longrightarrow C^\infty(M)$ be the divergence operator associated with $\mathbf{g}$. Its distributional extension is given by \[
(\operatorname{div}\mathbf{T})(\boldsymbol{\phi})
=
\mathbf{T}\bigl((\operatorname{div})'\boldsymbol{\phi}\bigr),
\qquad
\boldsymbol{\phi}\in\mathcal{D}(M).
\] In this case, $(\operatorname{div})'\boldsymbol{\phi}=-d\boldsymbol{\phi}$. Thus, if $W=\mathbb{K}$ and $\mathbf{T}=\mathbf{T}_{\mathbf{X}}$ for a smooth vector field $\mathbf{X}\in\mathfrak{X}(M)$, then \[
(\operatorname{div}\mathbf{T}_{\mathbf{X}})(\boldsymbol{\phi})
=
\int_M \boldsymbol{\phi}\,\operatorname{div}\mathbf{X}\,d\lambda_{\mathbf{g}}.
\] Using the Riemannian metric, \[
\int_M \boldsymbol{\phi}\,\operatorname{div}\mathbf{X}\,d\lambda_{\mathbf{g}}
=
-\int_M d\boldsymbol{\phi}(\mathbf{X})\,d\lambda_{\mathbf{g}}
=
-\int_M \langle \operatorname{grad}\boldsymbol{\phi},\mathbf{X}\rangle_{\mathbf{g}}\,d\lambda_{\mathbf{g}}.
\] Thus the bilinear formal transpose of divergence is written as $-d$, while the metric identifies it with the negative gradient.

\item Suppose that $\mathbf{E}$ is equipped with a bundle metric $\mathbf{h}_{\mathbf{E}}$, real or Hermitian according to the scalar field, and that $\nabla^{\mathbf{E}}$ is compatible with it. On $\nabla^2$, use the tensor connection induced by $\nabla^{\mathbf{E}}$ and the Levi--Civita connection of $\mathbf{g}$. Let $\Delta_B\colon \Gamma(\mathbf{E})\longrightarrow \Gamma(\mathbf{E})$ be the Bochner Laplacian associated with $\mathbf{g}$ and $\nabla^{\mathbf{E}}$, defined with the positive sign convention by $\Delta_B=(\nabla^{\mathbf{E}})_h^*\nabla^{\mathbf{E}}=-\operatorname{tr}_{\mathbf{g}}(\nabla^2)$. Its distributional extension is given by \[
(\Delta_B \mathbf{T})(\boldsymbol{\phi})
=
\mathbf{T}\bigl((\Delta_B)'\boldsymbol{\phi}\bigr),
\qquad
\boldsymbol{\phi}\in\mathcal{D}(M,\mathbf{E}^*).
\] If $W=\mathbb{K}$ and $\mathbf{T}=\mathbf{T}_{\mathbf{u}}$, then \[
(\Delta_B \mathbf{T}_{\mathbf{u}})(\boldsymbol{\phi})
=
\int_M \boldsymbol{\phi}(\Delta_B \mathbf{u})\,d\lambda_{\mathbf{g}}.
\] With this bundle metric $\mathbf{h}_{\mathbf{E}}$, define the test section $\boldsymbol{\phi}^{\sharp}\in\Gamma_c(\mathbf{E})$ by $\overline{\boldsymbol{\phi}^{\sharp}}:=\mathcal R_{\mathbf{E}}^{-1}\boldsymbol{\phi}$. Then \[
(\Delta_B \mathbf{T}_{\mathbf{u}})(\boldsymbol{\phi})
=
\int_M \mathbf{h}_{\mathbf{E}}(\Delta_B \mathbf{u},\boldsymbol{\phi}^{\sharp})\,d\lambda_{\mathbf{g}}
=
\int_M \mathbf{h}_{\mathbf{E}}(\mathbf{u},(\Delta_B)_h^*\boldsymbol{\phi}^{\sharp})\,d\lambda_{\mathbf{g}}.
\] Since $\Delta_B=(\nabla^{\mathbf{E}})_h^*\nabla^{\mathbf{E}}$, it is formally self-adjoint and satisfies $(\Delta_B)_h^*=\Delta_B$.

\item \emph{A Lorentzian digression, outside the Riemannian convention used throughout the rest of the chapter.} If $(M,\mathbf{g})$ is a closed Lorentzian manifold, the d'Alembert operator \[
\Box\colon C^\infty(M)\longrightarrow C^\infty(M),
\qquad
\Box:=\operatorname{tr}_{\mathbf{g}}(\nabla^2),
\] also induces an operator on distributions \[
\Box\colon \mathcal{D}'(M,W)\longrightarrow \mathcal{D}'(M,W),
\] defined by \[
(\Box \mathbf{T})(\boldsymbol{\phi})=\mathbf{T}(\Box'\boldsymbol{\phi}),
\qquad
\boldsymbol{\phi}\in\mathcal{D}(M).
\] In the regular scalar case, that is, if $W=\mathbb{K}$ and $\mathbf{T}=\mathbf{T}_{\mathbf{u}}$, we obtain \[
(\Box \mathbf{T}_{\mathbf{u}})(\boldsymbol{\phi})=\int_M \boldsymbol{\phi}\,\Box \mathbf{u}\,d\lambda_{\mathbf{g}}.
\]

In this case, $\lambda_{\mathbf{g}}$ denotes the positive Borel measure associated with the Lorentzian density: in coordinates, $d\lambda_{\mathbf{g}}=|\det(\mathbf{g})|^{\displaystyle\frac{1}{2}}\,dx$. No assertion is made that $\mathbf{g}$ defines a Riemannian distance, and none of the positive estimates specific to Riemannian geometry is used here. \item Finally, if $P\in \mathbf{PDO}^{(0)}(\mathbf{E},\mathbf{F})$ corresponds to a smooth section of $\operatorname{Hom}(\mathbf{E},\mathbf{F})$, then \[
(P\mathbf{T})(\boldsymbol{\phi})=\mathbf{T}(P'\boldsymbol{\phi}),
\qquad
\boldsymbol{\phi}\in\mathcal{D}(M,\mathbf{F}^*).
\] If $W=\mathbb{K}$ and $\mathbf{T}=\mathbf{T}_{\mathbf{u}}$, then \[
(P\mathbf{T}_{\mathbf{u}})(\boldsymbol{\phi})
=
\int_M \boldsymbol{\phi}(P\mathbf{u})\,d\lambda_{\mathbf{g}}.
\] In the presence of bundle metrics $\mathbf{h}_{\mathbf{E}}$ and $\mathbf{h}_{\mathbf{F}}$, define $\boldsymbol{\phi}^{\sharp}\in\Gamma_c(\mathbf{F})$ by $\overline{\boldsymbol{\phi}^{\sharp}}:=\mathcal R_{\mathbf{F}}^{-1}\boldsymbol{\phi}$. This can then be written as \[
(P\mathbf{T}_{\mathbf{u}})(\boldsymbol{\phi})
=
\int_M \mathbf{h}_{\mathbf{F}}(P\mathbf{u},\boldsymbol{\phi}^{\sharp})\,d\lambda_{\mathbf{g}}
=
\int_M \mathbf{h}_{\mathbf{E}}(\mathbf{u},P_h^*\boldsymbol{\phi}^{\sharp})\,d\lambda_{\mathbf{g}}.
\]

\end{enumerate}

In general, the action on distributions is determined by the bilinear formal transpose $P'$. For regular distributions with fixed bundle metrics, the preceding formulas can be rewritten using the formal Hermitian adjoint $P_h^*$, recovering the usual integration-by-parts identities. \end{example} \begin{example}[An example in which the space $W$ plays a role]\label{ej:distribuciones-en-haces-un-ejemplo-donde-el-espacio-si-interviene} Let $(M,\mathbf{g})$ be a Riemannian manifold with or without boundary, take the trivial bundle $\mathbf{E}=M\times\mathbb{K}$, and let $W=\mathbb{K}^2$. Fix a smooth function $u\in C^\infty(M)$ and a point $p\in M$. Define \[
T\colon \mathcal{D}(M)\longrightarrow \mathbb{K}^2
\] by \[
T(\phi)
:=
\left(
\int_M \phi u\,d\lambda_{\mathbf{g}},\,
\phi(p)
\right).
\] Then $T\in\mathcal{D}'(M,\mathbb{K}^2)$. The first component is a regular distribution, while the second is the Dirac delta at $p$.

If $P\in\mathbf{PDO}^{(m)}(M\times\mathbb{K},M\times\mathbb{K})$, the distributional extension of $P$ is given by \[
(PT)(\phi)=T(P'\phi).
\] Therefore \[
(PT)(\phi)
=
\left(
\int_M (P'\phi)u\,d\lambda_{\mathbf{g}},\,
(P'\phi)(p)
\right).
\] In particular, if $P=\Delta_{\mathbf{g}}$, whose interior formal transpose is itself, then \[
(\Delta_{\mathbf{g}} T)(\phi)
=
\left(
\int_M (\Delta_{\mathbf{g}}\phi)u\,d\lambda_{\mathbf{g}},\,
(\Delta_{\mathbf{g}}\phi)(p)
\right).
\] The first component represents the action of the Laplacian on the regular part, while the second represents its action on the Dirac delta. In this sense, taking $W=\mathbb{K}^2$ combines two different types of sources in a single distributional object: an extended source given by $u$ and a point source concentrated at $p$. \end{example} In the preceding examples, we mainly considered the scalar case $W=\mathbb{K}$. We now extend these constructions to distributions with values in a finite-dimensional vector space. \begin{definition}[Regular distribution with values in $W$]\label{def:distribuciones-en-haces-distribucion-regular-con-valores-en}\index{vector-valued regular distribution@vector-valued regular distribution} Let $(M,\mathbf{g})$ be a Riemannian manifold with or without boundary, let $\mathbf{E}\to M$ be a smooth vector bundle, and let $W$ be a finite-dimensional $\mathbb{K}$-vector space. A distribution \[
\mathbf{T}\colon \mathcal{D}(M,\mathbf{E}^*)\longrightarrow W
\] is called \textit{regular} if there exists $\mathbf{f}\in L^{1}_{\operatorname{loc}}(M,\mathbf{E}\otimes W)$ such that \[
\mathbf{T}(\boldsymbol{\phi})
=
\int_M (\boldsymbol{\phi}\otimes \operatorname{id}_W)(\mathbf{f})\,d\lambda_{\mathbf{g}},
\qquad \forall \boldsymbol{\phi}\in\mathcal{D}(M,\mathbf{E}^*).
\] \end{definition} \begin{remark}\label{obs:distribuciones-en-haces-caso-definicion-coincide-definicion-def-distribucion} When $W=\mathbb{K}$, this definition agrees with Definition~\ref{def: distribucion regular}. \end{remark} \begin{definition}\label{def:distribuciones-en-haces-espacio-vectorial-dimension-finita-distribucion-regular}\index{associated regular distribution@associated regular distribution} Let $(M,\mathbf{g})$ be a Riemannian manifold with or without boundary, let $\mathbf{E}\to M$ be a smooth vector bundle, let $W$ be a finite-dimensional $\mathbb{K}$-vector space, and let $\mathbf{u}\in L^1_{\operatorname{loc}}(M,\mathbf{E}\otimes W)$. Define the regular distribution associated with $\mathbf{u}$ by \[
\mathbf{T}_{\mathbf{u}}(\boldsymbol{\phi})
:=
\int_M (\boldsymbol{\phi}\otimes \operatorname{id}_W)(\mathbf{u})\,d\lambda_{\mathbf{g}},
\qquad
\boldsymbol{\phi}\in\mathcal{D}(M,\mathbf{E}^*).
\] \end{definition}

\begin{remark}\label{obs:distribuciones-en-haces-caso-seccion-suave-expresion-anterior-define} When $\mathbf{u}\in\Gamma(\mathbf{E}\otimes W)$ is a smooth section, the preceding expression defines a smooth regular distribution. In particular, this construction extends the classical association of smooth sections with regular distributions. \end{remark} \subsection{Local properties} A fundamental property of distributions is their local character. In particular, a distribution can be restricted to an open subset, reflecting the fact that its action is completely determined by the behavior of test functions on small neighborhoods. This construction is compatible with the vector bundle structure and will be used repeatedly below.

\begin{definition}[Restriction of distributions on bundles]\label{def:distribuciones-en-haces-restriccion-de-distribuciones-en-haces}\index{restriction of distributions on bundles@restriction of distributions on bundles} Let $M$ be a smooth manifold with or without boundary, $\mathbf{E}\longrightarrow M$ a smooth vector bundle, and $U\subseteq M$ an open subset. Let $\mathbf{T}\in \mathcal{D}'(M,\mathbf{E},W)$. Define the restriction of $\mathbf{T}$ to $U$ as the distribution \[
\mathbf{T}\restriction_U \in \mathcal{D}'(U,\mathbf{E},W)
\] given by \[
\mathbf{T}\restriction_U(\boldsymbol{\phi})
:=
\mathbf{T}(\widetilde{\boldsymbol{\phi}}),
\qquad
\boldsymbol{\phi}\in \mathcal{D}(U,\mathbf{E}^*),
\] where $\widetilde{\boldsymbol{\phi}}\in \mathcal{D}(M,\mathbf{E}^*)$ is the extension of $\boldsymbol{\phi}$ by zero outside $U$. \end{definition} \begin{lemma}[Well-definedness of restriction]\label{ejer:distribuciones-en-haces-variedad-suave-haz-vectorial-suave-espacio} Let $M$ be a smooth manifold with or without boundary, let $\mathbf{E}\longrightarrow M$ be a smooth vector bundle, let $W$ be a finite-dimensional $\mathbb{K}$-vector space, and let $U\subseteq M$ be open. If $\mathbf{T}\in\mathcal{D}'(M,\mathbf{E},W)$, the formula \[
\mathbf{T}\restriction_U(\boldsymbol{\phi}):=\mathbf{T}(\widetilde{\boldsymbol{\phi}}),
\qquad
\boldsymbol{\phi}\in\mathcal{D}(U,\mathbf{E}^*),
\] where $\widetilde{\boldsymbol{\phi}}$ denotes the extension of $\boldsymbol{\phi}$ by zero outside $U$, defines an element of $\mathcal{D}'(U,\mathbf{E},W)$. If $V\subseteq U$ is open, then \[
 (\mathbf{T}\restriction_U)\restriction_V=\mathbf{T}\restriction_V.
\]

\end{lemma}

\begin{proof} Denote extension by zero by \[
 e_U\colon\mathcal D(U,\mathbf{E}^*\restriction_U)
 \longrightarrow\mathcal D(M,\mathbf{E}^*)
\]. This map is linear. If $K\Subset U$, every $\boldsymbol{\phi}\in\Gamma_K(\mathbf{E}^*\restriction_U)$ vanishes on a neighborhood of the boundary of $U$, so $e_U\boldsymbol{\phi}$ is smooth and supported in $K$. For each $m\in\mathbb N_0$, its covariant derivatives agree with those of $\boldsymbol{\phi}$ on $U$ and vanish outside $U$; in particular, \[
 \|e_U\boldsymbol{\phi}\|_{C^m(K;M)}=\|\boldsymbol{\phi}\|_{C^m(K;U)}.
\] Thus the restriction of $e_U$ to every fixed-support stage is continuous. The final property of the $LF$ topology, given by Theorem~\ref{teo: topologia-final-lc}, implies that $e_U$ is continuous on the full test spaces.

Therefore $\mathbf{T}\restriction_U=\mathbf{T}\circ e_U$ is linear and continuous, that is, it belongs to $\mathcal D'(U,\mathbf{E},W)$. Finally, if $V\subseteq U$, extending by zero first from $V$ to $U$ and then from $U$ to $M$ gives exactly the extension by zero from $V$ to $M$. Composing with $\mathbf{T}$ yields the restriction identity. \end{proof} Restriction makes the local nature of distributions precise: the open subsets on which a distribution vanishes determine its support. \begin{definition}[Support of a distribution]\label{def:distribuciones-en-haces-soporte-de-una-distribucion}\index{support of a distribution@support of a distribution} Let $M$ be a smooth manifold with or without boundary, let $\mathbf{E}\to M$ be a smooth vector bundle, let $W$ be a finite-dimensional vector space, and let $\mathbf{T}\in \mathcal{D}'(M,\mathbf{E},W)$. Define the support of $\mathbf{T}$ as the set \[
\supp(\mathbf{T}):=
M\setminus \bigcup \{U\subseteq M \text{ open} \mid \mathbf{T}\restriction_U = 0\}.
\] \end{definition} The following proposition gives an equivalent definition: \begin{proposition}\label{prop: equivalencia definiciones soporte distribucion} Let $M$ be a smooth manifold with or without boundary, let $\mathbf{E}\longrightarrow M$ be a smooth vector bundle, and let $\mathbf{T}\in\mathcal{D}'(M,\mathbf{E},W)$. Then \[
\supp(\mathbf T)=\left\{x\in M\ \middle|
\begin{array}{l}
\text{for every open neighborhood }U\text{ of }x,\ \text{there exists }\\
\boldsymbol\phi\in\mathcal D(M,\mathbf E^*)\text{ with }\supp(\boldsymbol\phi)\subseteq U\\
\text{ and }\mathbf T(\boldsymbol\phi)\neq0
\end{array}\right\}.
\] \end{proposition}

\begin{proof} Define \[
A:=\bigcup\{U\subseteq M \text{ open}\mid \mathbf{T}\restriction_U=0\}.
\] Observe that $A$ is open, so its complement $M\setminus A$ is closed. We prove that $M\setminus A$ agrees with the set described in the second definition.

Let $x\in M\setminus A$. Consider an open neighborhood $U$ of $x$. We cannot have $\mathbf{T}\restriction_U=0$, since this would give $U\subseteq A$ and, in particular, $x\in A$, a contradiction. Thus $\mathbf{T}\restriction_U\neq 0$. This means there exists a test section $\boldsymbol{\phi}\in\mathcal{D}(U,\mathbf{E}^*)$ such that \[
\mathbf{T}\restriction_U(\boldsymbol{\phi})\neq 0.
\] If $\widetilde{\boldsymbol{\phi}}\in\mathcal{D}(M,\mathbf{E}^*)$ denotes the extension of $\boldsymbol{\phi}$ by zero outside $U$, then $\supp(\widetilde{\boldsymbol{\phi}})\subseteq U$ and \[
\mathbf{T}(\widetilde{\boldsymbol{\phi}})=\mathbf{T}\restriction_U(\boldsymbol{\phi})\neq 0.
\] We conclude that every open neighborhood of $x$ contains the support of a test section on which $\mathbf{T}$ does not vanish, so $x$ belongs to the second set.

Conversely, suppose $x$ belongs to the second set. If $x\in A$, there is an open subset $U\subseteq M$ such that $x\in U$ and $\mathbf{T}\restriction_U=0$. By hypothesis, however, there exists a test section $\boldsymbol{\phi}\in\mathcal{D}(M,\mathbf{E}^*)$ with $\supp(\boldsymbol{\phi})\subseteq U$ such that \[
\mathbf{T}(\boldsymbol{\phi})\neq 0.
\] Since $\supp(\boldsymbol{\phi})\subseteq U$, the restriction $\boldsymbol{\phi}\restriction_U$ belongs to $\mathcal{D}(U,\mathbf{E}^*)$, and its extension by zero to $M$ agrees with $\boldsymbol{\phi}$. Therefore \[
\mathbf{T}(\boldsymbol{\phi})=\mathbf{T}\restriction_U(\boldsymbol{\phi}\restriction_U)=0,
\] contradicting the choice of $\boldsymbol{\phi}$. Consequently, $x\notin A$.

We have proved that the two sets agree, and hence that the definitions are equivalent. \end{proof} Support makes the local character of distributions precise. In particular, a distribution does not ``detect'' test sections supported in regions where it vanishes. This is expressed by the following fundamental result.

\begin{proposition}\label{prop: separacion soportes} Let $M$ be a smooth manifold with or without boundary, let $\mathbf{E}\to M$ be a smooth vector bundle, and let $W$ be a finite-dimensional vector space. Let $\mathbf{T}\in\mathcal{D}'(M,\mathbf{E},W)$ and let $\boldsymbol{\phi}\in\mathcal{D}(M,\mathbf{E}^*)$. If $\supp(\boldsymbol{\phi})\cap\supp(\mathbf{T})=\emptyset$, then $\mathbf{T}(\boldsymbol{\phi})=0$.

\end{proposition}

\begin{proof} Let $K:=\supp(\boldsymbol{\phi})$, which is compact. For each $x\in K$, since $x\notin\supp(\mathbf{T})$, there is an open neighborhood $U_x$ of $x$ such that $\mathbf{T}\restriction_{U_x}=0$. The family $\{U_x\}_{x\in K}$ is an open cover of $K$, so there exist open subsets $U_1,\dots,U_k$ still covering $K$.

Let $\{\chi_j\}_{j\in\{1,\dots,k\}}$ be a partition of unity subordinate to $\{U_j\}_{j\in\{1,\dots,k\}}$. For each $j\in\{1,\dots,k\}$, define $\boldsymbol{\phi}_j:=\chi_j\,\boldsymbol{\phi}$. Then $\boldsymbol{\phi}=\displaystyle\sum_{j=1}^k \boldsymbol{\phi}_j$ and $\supp(\boldsymbol{\phi}_j)\subseteq U_j$ for every $j\in\{1,\dots,k\}$.

Since $\mathbf{T}\restriction_{U_j}=0$ and $\supp(\boldsymbol{\phi}_j)\subseteq U_j$, it follows that $\mathbf{T}(\boldsymbol{\phi}_j)=0$ for every $j\in\{1,\dots,k\}$. Therefore $\mathbf{T}(\boldsymbol{\phi})=\displaystyle\sum_{j=1}^k \mathbf{T}(\boldsymbol{\phi}_j)=0$.

\end{proof}

\begin{remark}\label{obs:distribuciones-en-haces-condicion-no-sustituirse-condicion-anule-ejemplo} The condition $\supp(\phi)\cap\supp(T)=\emptyset$ cannot be replaced by the condition that $\phi$ vanish on $\supp(T)$. For example, if $M=\mathbb{R}$, $\mathbf{E}=M\times \mathbb{K}$, and $T\in\mathcal{D}'(\mathbb{R},M\times \mathbb{K},\mathbb{K})$ is given by $T(\phi)=\phi'(0)$, then $\supp(T)=\{0\}$, but we may have $\phi(0)=0$ while $T(\phi)\neq 0$. \end{remark}

The preceding property naturally extends the action of a distribution to certain sections without compact support.

\begin{proposition}\label{prop: extension soporte compacto} Let $M$ be a smooth manifold with or without boundary, let $\mathbf{E}\to M$ be a smooth vector bundle, and let $W$ be a finite-dimensional vector space. Let $\mathbf{T}\in\mathcal{D}'(M,\mathbf{E},W)$ and let $\boldsymbol{\phi}\in \Gamma (M,\mathbf{E}^*)$. If $\supp(\mathbf{T})\cap\supp(\boldsymbol{\phi})$ is compact, the quantity $\mathbf{T}(\boldsymbol{\phi})$ can be defined consistently. \end{proposition}

\begin{proof} Let $K:=\supp(\mathbf{T})\cap\supp(\boldsymbol{\phi})$, which is compact by hypothesis. Since $M$ is a smooth manifold, there exists a relatively compact open subset $V\subseteq M$ such that $K\subseteq V$. For each $x\in K$, there is a relatively compact chart $V_x$ containing $x$; compactness of $K$ gives a finite subcover $V_{x_1},\dots,V_{x_N}$, and taking $V:=\displaystyle\bigcup_{j=1}^N V_{x_j}$ gives an open subset with $\overline{V}$ compact and $K\subseteq V$.

Since $K\subseteq V$ and $K$ is compact, choose an open subset $U\subseteq M$ such that \[
K\subseteq U\subseteq \overline{U}\subseteq V.
\] Now take a function $\sigma\in\mathcal{D}(M,\mathbb{K})$ such that $\sigma\equiv 1$ on $\overline{U}$ and $\supp(\sigma)\subseteq V$. Define \[
\mathbf{T}(\boldsymbol{\phi}):=\mathbf{T}(\sigma\boldsymbol{\phi}).
\] Observe that $\sigma\boldsymbol{\phi}\in\mathcal{D}(M,\mathbf{E}^*)$, so the preceding expression makes sense.

To see that the definition is independent of the choice of $\sigma$, let $\sigma$ and $\sigma'$ be two functions with the preceding property. Then $\sigma-\sigma'$ vanishes on a neighborhood of $K$, and in particular \[
\supp\bigl((\sigma-\sigma')\boldsymbol{\phi}\bigr)\cap\supp(\mathbf{T})=\emptyset.
\] By Proposition~\ref{prop: separacion soportes}, it follows that \[
\mathbf{T}\bigl((\sigma-\sigma')\boldsymbol{\phi}\bigr)=0.
\] Therefore \[
\mathbf{T}(\sigma\boldsymbol{\phi})-\mathbf{T}(\sigma'\boldsymbol{\phi})=\mathbf{T}\bigl((\sigma-\sigma')\boldsymbol{\phi}\bigr)=0,
\] proving that the definition is independent of the choice of $\sigma$. \end{proof}

\begin{remark}\label{obs:distribuciones-en-haces-proposicion-anterior-muestra-accion-distribucion-extenderse} The preceding proposition shows that the action of a distribution extends to smooth sections without compact support, provided the intersection with the support of the distribution is compact. This reinforces the idea that the behavior of a distribution is completely determined by its support. \end{remark}

The localization used later to define function spaces requires coordinate representations of distributions. The following construction also fixes the volume factor; omitting it would change the components of a regular distribution.

\begin{proposition}[Local representation of a distribution] \label{prop:representacion-local-distribucion-coordenadas} \index{distribution!coordinate representation} Let $(M,\mathbf{g})$ be a Riemannian manifold with or without boundary, let $\mathbf{E}\to M$ be a vector bundle of rank $r$, let $\mathbf{T}\in\mathcal D'(M,\mathbf{E})$, and let \[
\phi\colon U\longrightarrow\Omega\subseteq\mathbb R^n,
\qquad
\mathbf{e}=(\mathbf{e}_1,\dots,\mathbf{e}_r)
\] be a chart and a local frame of $\mathbf{E}$, with dual frame $\mathbf{e}^*=(\mathbf{e}^1,\dots,\mathbf{e}^r)$. If \[
J_\phi(x):=
\sqrt{\det\bigl(\mathbf{g}(\phi^{-1}(x))\bigr)},
\qquad x\in\Omega,
\] define the coordinate component of $\mathbf{T}$ by \begin{equation}
\bigl\langle T^a_{\phi,\mathbf{e}},\psi\bigr\rangle
:=
(\mathbf{T}\restriction_U)
\left(
\left(\frac{\psi}{J_\phi}\circ\phi\right)\mathbf{e}^a
\right),
\qquad
\psi\in\mathcal D(\Omega),
\label{eq:componente-local-distribucion}
\end{equation} for $a\in\{1,\dots,r\}$. Then $T^a_{\phi,\mathbf{e}}\in\mathcal D'(\Omega)$. If $\mathbf{u}=\displaystyle\sum_{a=1}^r u^a \mathbf{e}_a\in L^1_{\operatorname{loc}}(U,\mathbf{E})$ and $\mathbf{T}=\mathbf{T}_{\mathbf{u}}$, we have \[
T^a_{\phi,\mathbf{e}}
=T_{u^a\circ\phi^{-1}}
\quad\text{on }\Omega,
\] where the right-hand side is the Euclidean distribution defined using Lebesgue measure.

In the scalar case, we omit the frame and write $(\mathbf{T}\restriction_U)\circ\phi^{-1}$ for this distribution. If $\kappa=\phi^{-1}\colon \Omega\longrightarrow U$ is the corresponding parametrization, we also write $(\mathbf{T}\restriction_U)\circ\kappa$; both expressions denote the same object.

If $h\in C_c^\infty(U)$, each component of $h\mathbf{T}$ has compact support contained in $\Omega$ and admits a canonical extension by zero to \(\mathbb R^n\). Specifically, if $K:=\operatorname{supp}((h\mathbf{T})^a_{\phi,\mathbf{e}})\Subset\Omega$ and $\chi\in C_c^\infty(\Omega)$ satisfies $\chi=1$ on a neighborhood of $K$, this extension is given by \begin{equation}
\langle\widetilde{(h\mathbf{T})^a_{\phi,\mathbf{e}}},\psi\rangle
:=
\langle(h\mathbf{T})^a_{\phi,\mathbf{e}},\chi(\psi\restriction_\Omega)\rangle,
\qquad \psi\in\mathcal D(\mathbb R^n),
\label{eq:extension-cero-componente-distribucion}
\end{equation} and is independent of $\chi$. \end{proposition}

\begin{proof} The map \[
\psi\longmapsto
\left(\frac{\psi}{J_\phi}\circ\phi\right)\mathbf{e}^a
\] is linear and continuous from $\mathcal D(\Omega)$ to $\mathcal D(U,\mathbf{E}^*\restriction_U)$: on each compact set, the Leibniz and chain rules control all its derivatives. Continuity of $\mathbf{T}\restriction_U$ proves that \eqref{eq:componente-local-distribucion} defines a Euclidean distribution. For a regular distribution, use the change-of-variables formula \[
\int_U \eta\,d\lambda_{\mathbf{g}}
=
\int_\Omega (\eta\circ\phi^{-1})(x)J_\phi(x)\,d\lambda_n(x),
\] valid for every integrable function $\eta$ supported in $U$. Thus \[
\int_U
\left(\frac{\psi}{J_\phi}\circ\phi\right)u^a\,d\lambda_{\mathbf{g}}
=
\int_\Omega\psi(x)u^a(\phi^{-1}(x))\,dx,
\] proving the second assertion.

The identity $\operatorname{supp}(h\mathbf{T})\subseteq\operatorname{supp}(h)\Subset U$ follows directly from the definition of the product of a distribution with a smooth function and Proposition~\ref{prop: separacion soportes}. Thus the set $K$ in the last assertion is compact in $\Omega$. If $\chi_1$ and $\chi_2$ equal one near $K$, then $(\chi_1-\chi_2)(\psi\restriction_\Omega)$ has support separated from $K$; Proposition~\ref{prop: separacion soportes} shows that the two formulas give the same value. Continuity of \eqref{eq:extension-cero-componente-distribucion} is immediate on each space of test functions supported in a fixed compact set. \end{proof}

Support describes where a distribution ``lives.'' We can refine this idea by distinguishing points where the distribution is smooth from those where it has singularities. \begin{definition}[Singular support]\label{def:distribuciones-en-haces-soporte-singular}\index{singular support} Let $M$ be a smooth manifold with or without boundary, let $\mathbf{E}\to M$ be a smooth vector bundle, let $W$ be a finite-dimensional vector space, and let $\mathbf{T}\in\mathcal{D}'(M,\mathbf{E},W)$. Define the singular support of $\mathbf{T}$ as the set \[
\operatorname{sing\,supp}(\mathbf{T})
:=
\{x\in M \mid \text{for every neighborhood } U \text{ of } x,\
\mathbf{T}\restriction_U \neq \mathbf{T}_{\mathbf{u}} \text{ for every } \mathbf{u}\in\Gamma(U,\mathbf{E}\otimes W)\}.
\] Equivalently, $x\notin \operatorname{sing\,supp}(\mathbf{T})$ if and only if there exist an open neighborhood $U$ of $x$ and a smooth section $\mathbf{u}\in\Gamma(U,\mathbf{E}\otimes W)$ such that $\mathbf{T}\restriction_U = \mathbf{T}_{\mathbf{u}}$.

\end{definition}

\begin{proposition}\label{prop:distribuciones-en-haces-subconjunto-cerrado} Let $M$ be a smooth manifold with or without boundary, let $\mathbf{E}\to M$ be a smooth vector bundle, let $W$ be a finite-dimensional vector space, and let $\mathbf{T}\in\mathcal{D}'(M,\mathbf{E},W)$. Then $\operatorname{sing\,supp}(\mathbf{T})$ is a closed subset of $M$, and \[
\operatorname{sing\,supp}(\mathbf{T})\subseteq \supp(\mathbf{T}).
\] \end{proposition}

\begin{proof} Let $x\notin \operatorname{sing\,supp}(\mathbf{T})$. Then there exist an open neighborhood $U$ of $x$ and a smooth section $\mathbf{u}\in\Gamma(U,\mathbf{E}\otimes W)$ such that $\mathbf{T}\restriction_U=\mathbf{T}_{\mathbf{u}}$. In particular, for every $y\in U$, $\mathbf{T}\restriction_U$ agrees with a distribution induced by a smooth section on a neighborhood of $y$, so $y\notin \operatorname{sing\,supp}(\mathbf{T})$. This shows that $U\subseteq M\setminus \operatorname{sing\,supp}(\mathbf{T})$, and hence $M\setminus \operatorname{sing\,supp}(\mathbf{T})$ is open. We conclude that $\operatorname{sing\,supp}(\mathbf{T})$ is closed.

For the second assertion, let $x\notin\supp(\mathbf{T})$. There is then an open neighborhood $U$ of $x$ such that $\mathbf{T}\restriction_U=0$. Since $0=\mathbf{T}_0$ with $0\in\Gamma(U,\mathbf{E}\otimes W)$, $\mathbf{T}\restriction_U$ agrees with a distribution induced by a smooth section. Thus $x\notin \operatorname{sing\,supp}(\mathbf{T})$, implying $\operatorname{sing\,supp}(\mathbf{T})\subseteq \supp(\mathbf{T})$. \end{proof} \begin{example}\label{ej:distribuciones-en-haces-distribucion-delta-dirac-recordemos-distribucion-valores} Let $M$ be a smooth manifold with or without boundary, let $\mathbf{E}\to M$ be a smooth vector bundle, and let $p\in M$. Suppose that $\mathbf{E}_p^*\neq \{0\}$. Fix an auxiliary bundle metric $\mathbf{h}_{\mathbf{E}}$ on $\mathbf{E}$ and its dual metric $\mathbf{h}_{\mathbf{E}}^*$ on $\mathbf{E}^*$. Consider \[
\delta_{p}\colon \mathcal{D}(M,\mathbf{E}^{*}) \longrightarrow \mathbf{E}^{*}_{p},
\qquad
\delta_{p}(\boldsymbol{\phi}):=\boldsymbol{\phi}(p),
\] the Dirac delta distribution at $p$. Recall that $\delta_p$ is a distribution with values in $\mathbf{E}_p^*$, that is, \[
\delta_p\in \mathcal{D}'(M,\mathbf{E},\mathbf{E}_p^*).
\]

We claim that \[
\supp(\delta_p)=\{p\}.
\] If $\dim M\geq1$, we also have $\operatorname{sing\,supp}(\delta_p)=\{p\}$; if $\dim M=0$, its singular support is empty.

Let $U\subseteq M$ be an open subset such that $p\notin U$. Then, for every test section $\boldsymbol{\phi}\in\mathcal{D}(M,\mathbf{E}^*)$ with $\supp(\boldsymbol{\phi})\subseteq U$, we have $\boldsymbol{\phi}(p)=0$, and hence $\delta_p(\boldsymbol{\phi})=0$. This shows that $\delta_p\restriction_U=0$, giving \[
\supp(\delta_p)\subseteq \{p\}.
\]

On the other hand, if $U$ is any open subset containing $p$, there is a test section $\boldsymbol{\phi}\in\mathcal{D}(M,\mathbf{E}^*)$ with $\supp(\boldsymbol{\phi})\subseteq U$ such that $\boldsymbol{\phi}(p)\neq 0$. It suffices to take a nonzero element of $\mathbf{E}_p^*$ and use Lemma~\ref{lema:seccion-valor-prescrito}. Thus $\delta_p(\boldsymbol{\phi})=\boldsymbol{\phi}(p)\neq 0$, implying $p\in\supp(\delta_p)$. We conclude that \[
\supp(\delta_p)=\{p\}.
\]

If $\dim M=0$, the manifold is discrete, and Riemannian measure is counting measure. Let $\mathbf{v}$ be the section of $\mathbf{E}\otimes\mathbf{E}_p^*$ that vanishes outside $p$ and whose value at $p$ is the canonical tensor $\displaystyle\sum_{a=1}^r e_a\otimes e^a$, where $\{e_a\}_{a=1}^r$ is any basis of $\mathbf{E}_p$ and $\{e^a\}_{a=1}^r$ its dual basis. This tensor is independent of the basis, and the section is smooth because every point is open. For every test section, \[
\mathbf{T}_{\mathbf v}(\boldsymbol\phi)
=\sum_{a=1}^r\boldsymbol\phi(p)(e_a)e^a
=\boldsymbol\phi(p)=\delta_p(\boldsymbol\phi).
\] Hence $\delta_p$ is regular and has empty singular support.

Now suppose that $d:=\dim M\geq1$. Let us show that $p\in\operatorname{sing\,supp}(\delta_p)$. Suppose, for a contradiction, that there exist an open subset $U$ with $p\in U$ and a smooth section $\mathbf{u}\in\Gamma(U,\mathbf{E}\otimes \mathbf{E}_p^*)$ such that \[
\delta_p\restriction_U=\mathbf{T}_{\mathbf{u}}.
\] Then, for every test section $\boldsymbol{\phi}\in\mathcal{D}(U,\mathbf{E}^*)$, we would have \[
\boldsymbol{\phi}(p)=\delta_p(\boldsymbol{\phi})=\mathbf{T}_{\mathbf{u}}(\boldsymbol{\phi})
=
\int_U(\boldsymbol{\phi}\otimes\operatorname{id}_{\mathbf{E}_p^*})(\mathbf{u})\,d\lambda_{\mathbf{g}}.
\]

Consider an open subset $V\subseteq M$ such that $p\in V$, $\overline V$ is compact, and $\overline V\subseteq U$. We can choose it so that $\mathbf{E}^*$ admits a local trivialization over $V$ and a chart $(V,\varphi)$ exists. By Lemma~\ref{lema:seccion-valor-prescrito}, choose $\boldsymbol{\eta}\in\Gamma_c(V,\mathbf{E}^*)$ such that $\boldsymbol{\eta}(p)\neq0$, and extend it by zero outside $V$.

Let $H=\mathbb R^d$ if the chart is interior and $H=\mathbb H^d$ if it is a boundary chart. Since $\varphi(V)$ is open in $H$, choose a decreasing sequence of radii $r_n>0$, with $r_n\to0^{+}$, so that $\overline{B_{\mathrm{euc}}(\varphi(p),r_1)}\cap H\subseteq\varphi(V)$. Restricting Euclidean cutoff functions centered at $\varphi(p)$ to $H$ and transferring them by the chart gives $\rho_n\in C_c^\infty(V)$ with \[
0\leq\rho_n\leq1,\qquad \rho_n(p)=1,\qquad
\supp(\rho_n)\subseteq
\varphi^{-1}\bigl(B_{\mathrm{euc}}(\varphi(p),r_n)\cap H\bigr).
\] The weight $\sqrt{\det(\mathbf g)}\circ\varphi^{-1}$ is bounded on the compact coordinate set of radius $r_1$. Consequently, there exists $C_0>0$ independent of $n$ such that \[
\lambda_{\mathbf g}(\supp\rho_n)
\leq C_0\lambda_d\bigl(B_{\mathrm{euc}}(\varphi(p),r_n)\bigr)
\longrightarrow0,
\] since $d\geq1$. All these supports lie in $V\subseteq U$.

Now define \[
\boldsymbol{\phi}_n:=\rho_n\,\boldsymbol{\eta}.
\] Then $\boldsymbol{\phi}_n\in\mathcal{D}(U,\mathbf{E}^*)$ and \[
\boldsymbol{\phi}_n(p)=\boldsymbol{\eta}(p)\neq 0
\] for every $n$.

On the other hand, the section $\mathbf{s}:=(\boldsymbol{\eta}\otimes\operatorname{id}_{\mathbf{E}_p^*})(\mathbf{u})$ is smooth on $V$ and supported in the compact set $\operatorname{supp}\boldsymbol\eta\Subset V$. It is therefore bounded on that compact set and zero outside it. Thus there exists a constant $C>0$ such that $|\mathbf{s}(q)|_{\mathbf{h}_{\mathbf{E}}^*(p)}\leq C$ for every $q\in V$.

Since $\boldsymbol{\phi}_n=\rho_n\boldsymbol{\eta}$ and $\supp(\rho_n)\subseteq V$, we have \[
\int_U(\boldsymbol{\phi}_n\otimes\operatorname{id}_{\mathbf{E}_p^*})(\mathbf{u})\,d\lambda_{\mathbf{g}}
=
\int_V \rho_n\,\mathbf{s}\,d\lambda_{\mathbf{g}}.
\] By the triangle inequality for integrals, \[
\left|
\int_V \rho_n\,\mathbf{s}\,d\lambda_{\mathbf{g}}
\right|_{\mathbf{h}_{\mathbf{E}}^*(p)}
\leq
\int_V |\rho_n|\,|\mathbf{s}|_{\mathbf{h}_{\mathbf{E}}^*(p)}\,d\lambda_{\mathbf{g}}.
\] Since $0\leq \rho_n\leq 1$ and $|\mathbf{s}|_{\mathbf{h}_{\mathbf{E}}^*(p)}\leq C$ on $\supp(\rho_n)$, we obtain \[
\left|
\int_U(\boldsymbol{\phi}_n\otimes\operatorname{id}_{\mathbf{E}_p^*})(\mathbf{u})\,d\lambda_{\mathbf{g}}
\right|_{\mathbf{h}_{\mathbf{E}}^*(p)}
\leq
C\,\lambda_{\mathbf{g}}(\supp(\rho_n)).
\] Recall that $\lambda_{\mathbf{g}}(\supp(\rho_n))\longrightarrow 0$; therefore

\[
\int_U(\boldsymbol{\phi}_n\otimes\operatorname{id}_{\mathbf{E}_p^*})(\mathbf{u})\,d\lambda_{\mathbf{g}}
\longrightarrow 0.
\] This is a contradiction, since these integrals equal $\boldsymbol{\phi}_n(p)=\boldsymbol{\eta}(p)\neq 0$. Thus no such smooth section $\mathbf{u}$ exists, and it follows that $p\in\operatorname{sing\,supp}(\delta_p)$.

\end{example} We now examine further properties of support, including its relation to restrictions and partial differential operators. \begin{proposition}[Compatibility of restriction with support]\label{prop: soporte restriccion distribucion}\index{compatibility of restriction with support@compatibility of restriction with support} Let $M$ be a smooth manifold with or without boundary, let $\mathbf{E}\to M$ be a smooth vector bundle, and let $W$ be a finite-dimensional vector space. Let $\mathbf{T}\in\mathcal{D}'(M,\mathbf{E},W)$ and let $U\subseteq M$ be open. Then \[
\supp(\mathbf{T}\restriction_U)=\supp(\mathbf{T})\cap U.
\] \end{proposition}

\begin{proof} First prove that $\supp(\mathbf{T}\restriction_U)\subseteq \supp(\mathbf{T})\cap U$. Let $x\in U$ be such that $x\notin\supp(\mathbf{T})$. There is then an open neighborhood $V\subseteq M$ of $x$ such that $\mathbf{T}\restriction_V=0$. Since $U\cap V$ is an open neighborhood of $x$ in $U$, we have \[
(\mathbf{T}\restriction_U)\restriction_{U\cap V}=\mathbf{T}\restriction_{U\cap V}=0.
\] Therefore $x\notin\supp(\mathbf{T}\restriction_U)$.

Conversely, let $x\in \supp(\mathbf{T})\cap U$. Consider an open neighborhood $V\subseteq U$ of $x$. Since $U$ is open in $M$, we may also view $V$ as an open subset of $M$. Since $x\in\supp(\mathbf{T})$, we cannot have $\mathbf{T}\restriction_V=0$. Thus there is a test section $\boldsymbol{\phi}\in\mathcal{D}(V,\mathbf{E}^*)$ such that $\mathbf{T}\restriction_V(\boldsymbol{\phi})\neq 0$. But $\mathbf{T}\restriction_V=(\mathbf{T}\restriction_U)\restriction_V$, so $(\mathbf{T}\restriction_U)\restriction_V\neq 0$. Since this holds for every open neighborhood $V$ of $x$ in $U$, it follows that $x\in\supp(\mathbf{T}\restriction_U)$.

We conclude that $\supp(\mathbf{T}\restriction_U)=\supp(\mathbf{T})\cap U$. \end{proof}

\begin{proposition}[Compatibility of restriction with differential operators] \label{prop: restriccion conmuta pdo distribuciones} \index{compatibility of restriction with differential operators@compatibility of restriction with differential operators} Let $M$ be a smooth manifold without boundary, let $\mathbf{E},\mathbf{F}\to M$ be smooth vector bundles, and let $W$ be a finite-dimensional vector space. Let $P\in\mathbf{PDO}^{(m)}(\mathbf{E},\mathbf{F})$, let $\mathbf{T}\in\mathcal D'(M,\mathbf{E},W)$, and let $U\subseteq M$ be open. If $P_U$ is the restriction of $P$ to $U$, then \begin{equation}
\label{eq:restriccion-pdo-distribuciones}
 (P\mathbf{T})|_U=P_U(\mathbf{T}|_U)
 \qquad\text{on }\mathcal D'(U,\mathbf{F}|_U,W).
\end{equation} \end{proposition}

\begin{proof} Compatibility of the formal transpose with restrictions gives $(P_U)'=(P')_U$. Let $\boldsymbol{\phi}\in\mathcal D(U,\mathbf{F}^*|_U)$ and denote its extension by zero by $\widetilde{\boldsymbol{\phi}}$; it is smooth because $\operatorname{supp}\boldsymbol{\phi}\Subset U$. Locality of $P'$ implies that $P'\widetilde{\boldsymbol{\phi}}$ is the extension by zero of $(P_U)'\boldsymbol{\phi}$. Therefore \begin{align*}
 ((P\mathbf{T})|_U)(\boldsymbol{\phi})
 &=(P\mathbf{T})(\widetilde{\boldsymbol{\phi}})
 =\mathbf{T}(P'\widetilde{\boldsymbol{\phi}})\\
 &=(\mathbf{T}|_U)((P_U)'\boldsymbol{\phi})
 =(P_U(\mathbf{T}|_U))(\boldsymbol{\phi}),
\end{align*} which is \eqref{eq:restriccion-pdo-distribuciones}. \end{proof}

\begin{proposition}[Locality of differential operators on distributions]\label{prop: soporte pdo distribuciones}\index{locality of differential operators on distributions} Let $M$ be a smooth manifold without boundary, let $\mathbf{E},\mathbf{F}\to M$ be smooth vector bundles, and let $W$ be a finite-dimensional vector space. Let $P\in\mathbf{PDO}^{(m)}(\mathbf{E},\mathbf{F})$ and let $\mathbf{T}\in\mathcal{D}'(M,\mathbf{E},W)$. Then \[
\supp(P\mathbf{T})\subseteq \supp(\mathbf{T}).
\] \end{proposition}

\begin{proof} Let $x\notin\supp(\mathbf{T})$. There is then an open neighborhood $U$ of $x$ such that $\mathbf{T}\restriction_U=0$.

Let $\boldsymbol{\phi}\in\mathcal{D}(U,\mathbf{F}^*)$ and let $\widetilde{\boldsymbol{\phi}}\in\mathcal{D}(M,\mathbf{F}^*)$ be its extension by zero outside $U$. By Proposition~\ref{prop: restriccion conmuta pdo distribuciones}, we have \[
(P\mathbf{T})\restriction_U(\boldsymbol{\phi})
=
\mathbf{T}\restriction_U\bigl((P'\widetilde{\boldsymbol{\phi}})\restriction_U\bigr).
\] Since $\mathbf{T}\restriction_U=0$, we obtain \[
(P\mathbf{T})\restriction_U(\boldsymbol{\phi})=0.
\] This holds for every $\boldsymbol{\phi}\in\mathcal{D}(U,\mathbf{F}^*)$, so $(P\mathbf{T})\restriction_U=0$. Consequently, $x\notin\supp(P\mathbf{T})$.

Therefore \[
\supp(P\mathbf{T})\subseteq \supp(\mathbf{T}).
\] \end{proof} \begin{lemma}[Restriction of regular distributions]\label{lem: restriccion distribucion regular}\index{restriction of regular distributions@restriction of regular distributions} Let $(M,\mathbf{g})$ be a Riemannian manifold with or without boundary, let $\mathbf{E}\to M$ be a smooth vector bundle, and let $W$ be a finite-dimensional vector space. Let $\mathbf{u}\in\Gamma(\mathbf{E}\otimes W)$ and let $U\subseteq M$ be open. Then \[
(\mathbf{T}_{\mathbf{u}})\restriction_U=\mathbf{T}_{\mathbf{u}\restriction_U}
\] as elements of $\mathcal{D}'(U,\mathbf{E}\restriction_U,W)$. \end{lemma}

\begin{proof} Let $\boldsymbol{\phi}\in\mathcal{D}(U,\mathbf{E}^*\restriction_U)$. By the definition of restriction of a distribution, \[
(\mathbf{T}_{\mathbf{u}})\restriction_U(\boldsymbol{\phi})=\mathbf{T}_{\mathbf{u}}(\widetilde{\boldsymbol{\phi}}),
\] where $\widetilde{\boldsymbol{\phi}}\in\mathcal{D}(M,\mathbf{E}^*)$ denotes the extension of $\boldsymbol{\phi}$ by zero outside $U$.

By the definition of $\mathbf{T}_{\mathbf{u}}$, we have \[
\mathbf{T}_{\mathbf{u}}(\widetilde{\boldsymbol{\phi}})
=
\int_M(\widetilde{\boldsymbol{\phi}}\otimes \operatorname{id}_W)(\mathbf{u})\,d\lambda_{\mathbf{g}}.
\] Since $\widetilde{\boldsymbol{\phi}}$ vanishes outside $U$, this integral reduces to \[
\int_U(\boldsymbol{\phi}\otimes \operatorname{id}_W)(\mathbf{u}\restriction_U)\,d\lambda_{\mathbf{g}}.
\] But the last expression is precisely \[
\mathbf{T}_{\mathbf{u}\restriction_U}(\boldsymbol{\phi}).
\] Therefore $(\mathbf{T}_{\mathbf{u}})\restriction_U(\boldsymbol{\phi})=\mathbf{T}_{\mathbf{u}\restriction_U}(\boldsymbol{\phi})$ for every $\boldsymbol{\phi}\in\mathcal{D}(U,\mathbf{E}^*\restriction_U)$, proving the equality. \end{proof} \begin{proposition}[Smooth regular distributions and support]\label{prop: soporte distribucion regular suave}\index{smooth regular distributions and support} Let $(M,\mathbf{g})$ be a Riemannian manifold with or without boundary, let $\mathbf{E}\to M$ be a smooth vector bundle, and let $W$ be a finite-dimensional vector space. Let $\mathbf{u}\in\Gamma(\mathbf{E}\otimes W)$ and let $\mathbf{T}_{\mathbf{u}}\in\mathcal{D}'(M,\mathbf{E},W)$ be the smooth regular distribution induced by $\mathbf{u}$. Then \[
\supp(\mathbf{T}_{\mathbf{u}})=\supp(\mathbf{u}).
\] \end{proposition}

\begin{proof} Let $x\notin\supp(\mathbf{u})$. There is then an open neighborhood $U$ of $x$ such that $\mathbf{u}\restriction_U=0$. If $\boldsymbol{\phi}\in\mathcal{D}(U,\mathbf{E}^*)$, we have \[
(\mathbf{T}_{\mathbf{u}}\restriction_U)(\boldsymbol{\phi})
=
\mathbf{T}_{\mathbf{u}\restriction_U}(\boldsymbol{\phi})
=
\int_U(\boldsymbol{\phi}\otimes\operatorname{id}_W)(\mathbf{u}\restriction_U)\,d\lambda_{\mathbf{g}}
=
0.
\] It follows that $\mathbf{T}_{\mathbf{u}}\restriction_U=0$, and hence $x\notin\supp(\mathbf{T}_{\mathbf{u}})$. This proves that $\supp(\mathbf{T}_{\mathbf{u}})\subseteq\supp(\mathbf{u})$.

Conversely, let $x\in\supp(\mathbf{u})$. We show that $x\in\supp(\mathbf{T}_{\mathbf{u}})$. Let $U$ be an arbitrary open neighborhood of $x$. Since $x\in\supp(\mathbf{u})$, there exists $y\in U$ such that $\mathbf{u}(y)\neq 0$.

Since $\mathbf{u}(y)\in \mathbf{E}_y\otimes W$ is nonzero, there exist $\ell_0\in \mathbf{E}_y^*$ and $\omega\in W^*$ such that \[
\omega\bigl((\ell_0\otimes \operatorname{id}_W)(\mathbf{u}(y))\bigr)\neq 0.
\] Applying Lemma~\ref{lema:seccion-valor-prescrito} to the bundle $\mathbf{E}^*\to M$, choose a smooth section $\boldsymbol{\ell}\in\Gamma(\mathbf{E}^*)$ with $\supp(\boldsymbol{\ell})\subseteq U$ and $\boldsymbol{\ell}(y)=\ell_0$.

Multiplying $\boldsymbol{\ell}$ by a nonzero scalar if necessary, we may assume that \[
\operatorname{Re}\,
\omega\bigl((\boldsymbol{\ell}(y)\otimes \operatorname{id}_W)(\mathbf{u}(y))\bigr)>0,
\] where, if $\mathbb{K}=\mathbb{R}$, $\operatorname{Re}$ is understood as the identity. By continuity, there is a relatively compact open subset $V\subseteq U$, with $y\in V$, such that \[
\operatorname{Re}\,
\omega\bigl((\boldsymbol{\ell}(z)\otimes \operatorname{id}_W)(\mathbf{u}(z))\bigr)>0
\qquad
\forall z\in V.
\]

Now let $\rho\in\mathcal{D}(U,\mathbb{K})$ be a cutoff function such that $\rho\geq 0$, $\rho\not\equiv 0$, and $\supp(\rho)\subseteq V$. Define \[
\boldsymbol{\phi}:=\rho\,\boldsymbol{\ell}\restriction_U\in\mathcal{D}(U,\mathbf{E}^*).
\] Then \[
\omega\bigl(\mathbf{T}_{\mathbf{u}}\restriction_U(\boldsymbol{\phi})\bigr)
=
\omega\left(
\int_U(\boldsymbol{\phi}\otimes\operatorname{id}_W)(\mathbf{u})\,d\lambda_{\mathbf{g}}
\right).
\] Since $W$ is finite dimensional, we may move $\omega$ inside the integral, obtaining \[
\omega\bigl(\mathbf{T}_{\mathbf{u}}\restriction_U(\boldsymbol{\phi})\bigr)
=
\int_U
\omega\bigl((\boldsymbol{\phi}\otimes\operatorname{id}_W)(\mathbf{u})\bigr)\,d\lambda_{\mathbf{g}}.
\] Using the definition of $\boldsymbol{\phi}$ and the fact that $\supp(\rho)\subseteq V$, this reduces to \[
\omega\bigl(\mathbf{T}_{\mathbf{u}}\restriction_U(\boldsymbol{\phi})\bigr)
=
\int_V
\rho\,
\omega\bigl((\boldsymbol{\ell}\otimes\operatorname{id}_W)(\mathbf{u})\bigr)\,d\lambda_{\mathbf{g}}.
\] Taking real parts, \[
\operatorname{Re}\,\omega\bigl(\mathbf{T}_{\mathbf{u}}\restriction_U(\boldsymbol{\phi})\bigr)
=
\int_V
\rho\,
\operatorname{Re}\,
\omega\bigl((\boldsymbol{\ell}\otimes\operatorname{id}_W)(\mathbf{u})\bigr)\,d\lambda_{\mathbf{g}}.
\] The integrand is nonnegative and does not vanish identically, since $\rho\geq 0$, $\rho\not\equiv 0$, and \[
\operatorname{Re}\,
\omega\bigl((\boldsymbol{\ell}(z)\otimes\operatorname{id}_W)(\mathbf{u}(z))\bigr)>0
\qquad
\forall z\in V.
\] Consequently, \[
\operatorname{Re}\,\omega\bigl(\mathbf{T}_{\mathbf{u}}\restriction_U(\boldsymbol{\phi})\bigr)>0,
\] and in particular $\mathbf{T}_{\mathbf{u}}\restriction_U(\boldsymbol{\phi})\neq 0$. Therefore $\mathbf{T}_{\mathbf{u}}\restriction_U\neq 0$.

Since $U$ was an arbitrary neighborhood of $x$, we conclude that $x\in\supp(\mathbf{T}_{\mathbf{u}})$. Consequently, \[
\supp(\mathbf{T}_{\mathbf{u}})=\supp(\mathbf{u}).
\] \end{proof}

\begin{proposition}[Compatibility of restriction with singular support]\label{prop: soporte singular restriccion}\index{compatibility of restriction with singular support@compatibility of restriction with singular support} Let $M$ be a smooth manifold with or without boundary, let $\mathbf{E}\to M$ be a smooth vector bundle, and let $W$ be a finite-dimensional vector space. Let $\mathbf{T}\in\mathcal{D}'(M,\mathbf{E},W)$ and let $U\subseteq M$ be open. Then \[
\operatorname{sing\,supp}(\mathbf{T}\restriction_U)=\operatorname{sing\,supp}(\mathbf{T})\cap U.
\] \end{proposition}

\begin{proof} Let $x\in U$ be such that $x\notin\operatorname{sing\,supp}(\mathbf{T})$. There then exist an open neighborhood $V\subseteq M$ of $x$ and a smooth section $\mathbf{u}\in\Gamma(V,\mathbf{E}\otimes W)$ such that \[
\mathbf{T}\restriction_V=\mathbf{T}_{\mathbf{u}}.
\] Restricting to $U\cap V$ gives \[
(\mathbf{T}\restriction_U)\restriction_{U\cap V}=\mathbf{T}\restriction_{U\cap V}=(\mathbf{T}_{\mathbf{u}})\restriction_{U\cap V}=\mathbf{T}_{\mathbf{u}\restriction_{U\cap V}}.
\] Therefore $x\notin\operatorname{sing\,supp}(\mathbf{T}\restriction_U)$. This proves that \[
\operatorname{sing\,supp}(\mathbf{T}\restriction_U)\subseteq \operatorname{sing\,supp}(\mathbf{T})\cap U.
\]

Conversely, let $x\in U$ be such that $x\notin\operatorname{sing\,supp}(\mathbf{T}\restriction_U)$. There then exist an open neighborhood $V\subseteq U$ of $x$ and a smooth section $\mathbf{u}\in\Gamma(V,\mathbf{E}\otimes W)$ such that \[
(\mathbf{T}\restriction_U)\restriction_V=\mathbf{T}_{\mathbf{u}}.
\] Since $V$ is also open in $M$, we have $\mathbf{T}\restriction_V=(\mathbf{T}\restriction_U)\restriction_V=\mathbf{T}_{\mathbf{u}}$. Therefore $x\notin\operatorname{sing\,supp}(\mathbf{T})$. This proves the reverse containment.

We conclude that \[
\operatorname{sing\,supp}(\mathbf{T}\restriction_U)=\operatorname{sing\,supp}(\mathbf{T})\cap U.
\] \end{proof}

\begin{proposition}[Locality of differential operators and singular support]\label{prop: soporte singular pdo}\index{locality of differential operators and singular support} Let $M$ be a smooth manifold without boundary, let $\mathbf{E},\mathbf{F}\to M$ be smooth vector bundles, and let $W$ be a finite-dimensional vector space. Let $P\in\mathbf{PDO}^{(m)}(\mathbf{E},\mathbf{F})$ and let $\mathbf{T}\in\mathcal{D}'(M,\mathbf{E},W)$. Then \[
\operatorname{sing\,supp}(P\mathbf{T})\subseteq \operatorname{sing\,supp}(\mathbf{T}).
\] \end{proposition}

\begin{proof} Let $x\notin\operatorname{sing\,supp}(\mathbf{T})$. There then exist an open neighborhood $U$ of $x$ and a smooth section $\mathbf{u}\in\Gamma(U,\mathbf{E}\otimes W)$ such that \[
\mathbf{T}\restriction_U=\mathbf{T}_{\mathbf{u}}.
\] We show that $(P\mathbf{T})\restriction_U$ is also a smooth regular distribution.

Let $\boldsymbol{\phi}\in\mathcal{D}(U,\mathbf{F}^*)$ and let $\widetilde{\boldsymbol{\phi}}\in\mathcal{D}(M,\mathbf{F}^*)$ be its extension by zero outside $U$. By Proposition~\ref{prop: restriccion conmuta pdo distribuciones}, we have \[
(P\mathbf{T})\restriction_U(\boldsymbol{\phi})
=
\mathbf{T}\restriction_U\bigl((P'\widetilde{\boldsymbol{\phi}})\restriction_U\bigr).
\] Using $\mathbf{T}\restriction_U=\mathbf{T}_{\mathbf{u}}$, we obtain \[
(P\mathbf{T})\restriction_U(\boldsymbol{\phi})
=
\mathbf{T}_{\mathbf{u}}\bigl((P'\widetilde{\boldsymbol{\phi}})\restriction_U\bigr)
=
\int_U
\bigl((P'\widetilde{\boldsymbol{\phi}})\restriction_U\otimes\operatorname{id}_{W}\bigr)(\mathbf{u})\,d\lambda_{\mathbf{g}}.
\]

Now observe that $\widetilde{\boldsymbol{\phi}}$ agrees with $\boldsymbol{\phi}$ on $U$, so $\widetilde{\boldsymbol{\phi}}-\boldsymbol{\phi}=0$ on $U$. Since $P'$ is a differential operator, it is local, and therefore \[
P'(\widetilde{\boldsymbol{\phi}}-\boldsymbol{\phi})=0
\qquad \text{on } U.
\] Consequently, \[
(P'\widetilde{\boldsymbol{\phi}})\restriction_U=P'\boldsymbol{\phi},
\] where $P'\boldsymbol{\phi}$ is computed in $U$.

Substitute this into the preceding integral. To justify the resulting identity, work locally: write \[
\mathbf{u}=\displaystyle\sum_{a=1}^{N}\mathbf{u}_a\otimes w_a,
\] where $\{w_1,\dots,w_N\}$ is a basis of $W$ and each $\mathbf{u}_a$ is a smooth local section of $\mathbf{E}$. Then \[
\bigl((P'\widetilde{\boldsymbol{\phi}})\restriction_U\otimes\operatorname{id}_{W}\bigr)(\mathbf{u})
=
\displaystyle\sum_{a=1}^{N}(P'\boldsymbol{\phi})(\mathbf{u}_a)\,w_a.
\] Applying the definition of the bilinear formal transpose (on $U$), we have \[
\int_U (P'\boldsymbol{\phi})(\mathbf{u}_a)\,d\lambda_{\mathbf{g}}
=
\int_U \boldsymbol{\phi}(P\mathbf{u}_a)\,d\lambda_{\mathbf{g}}
\qquad
\text{for every }a\in\{1,\dots,N\}.
\] Therefore \[
\int_U
\bigl((P'\widetilde{\boldsymbol{\phi}})\restriction_U\otimes\operatorname{id}_{W}\bigr)(\mathbf{u})\,d\lambda_{\mathbf{g}}
=
\displaystyle\sum_{a=1}^{N}
\left(
\int_U \boldsymbol{\phi}(P\mathbf{u}_a)\,d\lambda_{\mathbf{g}}
\right)w_a.
\]

On the other hand, since $P$ acts only on the component in $\mathbf{E}$, we have \[
(P\otimes\operatorname{id}_{W})\mathbf{u}
=
\displaystyle\sum_{a=1}^{N}P\mathbf{u}_a\otimes w_a
\in\Gamma(U,\mathbf{F}\otimes W).
\] It follows that \[
\displaystyle\sum_{a=1}^{N}
\left(
\int_U \boldsymbol{\phi}(P\mathbf{u}_a)\,d\lambda_{\mathbf{g}}
\right)w_a
=
\int_U
\bigl(\boldsymbol{\phi}\otimes\operatorname{id}_{W}\bigr)
\bigl((P\otimes\operatorname{id}_{W})\mathbf{u}\bigr)\,d\lambda_{\mathbf{g}}.
\]

Consequently, \[
(P\mathbf{T})\restriction_U(\boldsymbol{\phi})
=
\int_U
\bigl(\boldsymbol{\phi}\otimes\operatorname{id}_{W}\bigr)
\bigl((P\otimes\operatorname{id}_{W})\mathbf{u}\bigr)\,d\lambda_{\mathbf{g}}
=
\mathbf{T}_{(P\otimes\operatorname{id}_{W})\mathbf{u}}(\boldsymbol{\phi}).
\]

Since this holds for every $\boldsymbol{\phi}\in\mathcal{D}(U,\mathbf{F}^*)$, we conclude that \[
(P\mathbf{T})\restriction_U
=
\mathbf{T}_{(P\otimes\operatorname{id}_{W})\mathbf{u}}.
\] In particular, $(P\mathbf{T})\restriction_U$ is a smooth regular distribution on a neighborhood of $x$, implying \[
x\notin\operatorname{sing\,supp}(P\mathbf{T}).
\] Since $x$ was arbitrary, we obtain \[
\operatorname{sing\,supp}(P\mathbf{T})\subseteq \operatorname{sing\,supp}(\mathbf{T}).
\] \end{proof}

\subsection{External products and parameter-dependent sections} To study operators from a more global viewpoint, it is useful to consider sections depending on two variables. In the Euclidean case, these arise when operators are written using functions of two variables, usually called kernels. For manifolds and vector bundles, the corresponding expression must be formulated intrinsically.

\begin{remark}[Products with boundary and corners]
\label{obs:producto-variedades-con-esquinas-prueba}
Theorem~\ref{teo:producto-variedades-con-esquinas} constructs the
product for factors with arbitrary corners, and
Corollary~\ref{cor:producto-dos-variedades-frontera} specifies the case
of two factors with smooth boundary. We use its product structure,
with smoothness defined by local Euclidean extensions.
Lemma~\ref{lem:esquinas-jets-extension-local} identifies this condition
with continuous extension to the faces of all interior derivatives.
Test functions may be nonzero on faces and their intersections.
This convention retains distributions supported there. We next
justify the topologies needed to use these products in the kernel
theorem.
\end{remark}

\begin{proposition}[Test spaces on manifolds with corners]
\label{prop:prueba-esquinas-Frechet-LF}
Let $M$ be a smooth manifold with corners, let $\mathbf E\to M$ be a
vector bundle of finite rank, and let $K\subseteq M$ be compact.
The $C^j(K)$ norms, defined by smooth metrics and connections on the
ordinary tangent bundle, make $\Gamma_K(\mathbf E)$ a Fréchet space
whose topology is independent of these auxiliary data.
For $k<\infty$, $\Gamma_K^k(\mathbf E)$ is Banach. If $K\subseteq L$
are compact, the fixed-support inclusions are embeddings with closed
image. An exhaustion $K_j\subseteq\operatorname{Int}_M K_{j+1}$ defines
\[
 \mathcal D(M,\mathbf E):=\varinjlim_j\Gamma_{K_j}(\mathbf E)
\]
as a complete strict $LF$ space, independent of the exhaustion.
Moreover, if $K\subseteq\operatorname{Int}_M L$, every section in
$\Gamma_K^k(\mathbf E)$ is approximated in the $C^k$ norm by smooth
sections supported in $L$.
\end{proposition}

\begin{proof}
Cover $K$ by finitely many charts and frames, and choose cutoffs
$\psi_i$ compactly supported in their domains with sum one near $K$,
by Proposition~\ref{prop:esquinas-particiones-unidad}. The coordinate
components $f_i^a$ of $\psi_i\mathbf u$ extend by zero to $G_{n,r_i}$
and have supports in fixed compact sets $C_i$. The triangular
identities between ordinary and covariant derivatives in
Proposition~\ref{prop:sobolev-esquinas-localizacion-densidad}, now in
uniform norm, give positive constants for each $j$ with
\[
 c_j\|\mathbf u\|_{C^j(K)}
 \leq\sum_{i,a}\|f_i^a\|_{C^j(G_{n,r_i})}
 \leq C_j\|\mathbf u\|_{C^j(K)}.
\]
The constants only use coefficient bounds on compact sets. Thus
they also prove independence of metrics and connections.

Let $(\mathbf u_b)$ be Cauchy for all norms. The operator $E_{n,r_i}$
in Lemma~\ref{lem:extension-suave-semiespacio-soporte-fijo} extends
$f_{i,b}^a$ to smooth Euclidean functions supported in a fixed compact
set $\widehat C_i$ and has bounds in each $C^j$ norm. The extensions
are therefore Cauchy in $\mathcal D_{\widehat C_i}(\mathbb R^n)$,
complete by Corollary~\ref{cor:D-K-euclidiano-nuclear}. Their limits
are smooth Euclidean functions. Restriction gives smooth limits
$f_i^a$ supported in $C_i$, since uniform convergence preserves
vanishing off that compact set. Reconstruction and frame-change
identities pass to the limit in every derivative. Reconstructing the
components in their frames and extending by zero gives a smooth
section $\mathbf u$ supported in $K$, and the preceding bound gives
$\mathbf u_b\to\mathbf u$ in each $C^j(K)$ norm.
The countable family of norms separates points; metrization by
seminorms in Proposition~\ref{prop: metrizacion por seminormas} gives
a complete metric. This proves the Fréchet assertion.

For finite order $k$, use the finite $C^k$ extension in
Lemma~\ref{lem:holder-esquinas-independencia} instead of the smooth
extension. A Cauchy sequence has $C^k$ limits in each component,
reconstructed in the same way; this proves completeness of
$\Gamma_K^k$. These extensions, followed by convolution, approximate
the localized components in the $C^k$ norm. If
$K\subseteq\operatorname{Int}_M L$, choose all localization and
reconstruction cutoffs supported in $\operatorname{Int}_M L$.
Reconstructing the smooth approximations gives sections supported in
$L$ converging to the original section. Smooth density with exactly
the same support $K$ is not asserted.

If $K\subseteq L$, a section supported in $K$ and all its derivatives
vanish off $K$. Hence
$\|\mathbf u\|_{C^j(L)}=\|\mathbf u\|_{C^j(K)}$, proving that the
inclusion induces the same topology. If a sequence converges in the
$L$ stage, its uniform limit vanishes off $K$, so its support stays
in $K$. The image is closed for both finite and infinite order.
The required exhaustion exists by
Proposition~\ref{prop:esquinas-particiones-unidad}. We thus have a
countable inductive limit of Fréchet spaces with closed inclusions
inducing their topologies; it is strict $LF$ and is complete by
Theorem~\ref{teo: propiedades LF} of Appendix~B.
Two exhaustions are cofinal: each compact set of one lies in a compact
set of the other, since the increasing stage interiors cover $M$.
Continuous inclusions between stages give continuity of the identity
in both directions for the final topologies. This proves independence
of the exhaustion.
\end{proof}

For a manifold with corners define $\mathcal D'(M,\mathbf E,W)$ as
the continuous linear maps from $\mathcal D(M,\mathbf E^*)$ to $W$,
with pointwise convergence topology, exactly as in
Definition~\ref{def:distribuciones-en-haces-distribuciones}.
Continuity is equivalent to continuity on each fixed-support stage.
On a stage, a continuous functional has a bound
$|\mathbf T(\boldsymbol\varphi)|\leq C\|\boldsymbol\varphi\|_{C^m(K)}$
for some finite order $m$, because the norms form an increasing
family defining its topology. Evaluation at a point of a face or
corner is continuous already for the order-zero norm. In the
tensorization and kernel theorems below, factors may be manifolds
with corners of any depth.

If $\mathbf{E}\longrightarrow M$ and $\mathbf{F}\longrightarrow N$ are smooth vector bundles, a section depending on two variables and taking values in $\mathbf{E}_x\otimes \mathbf{F}_y^*$ for each $(x,y)\in M\times N$ must be understood as a section of a vector bundle over $M\times N$.

Pullback by the canonical projections of the product gives vector bundles over $M\times N$ from bundles defined over each factor. This allows us to introduce the external product.

\begin{definition}[External product of vector bundles]\label{def:distribuciones-en-haces-producto-externo-de-haces-vectoriales}\index{external product of vector bundles} Let $M$ and $N$ be smooth manifolds with corners, and let $\mathbf{E}\longrightarrow M$ and $\mathbf{F}\longrightarrow N$ be smooth vector bundles. Denote the canonical projections by $\pi_M\colon M\times N\longrightarrow M$ and $\pi_N\colon M\times N\longrightarrow N$.

Define the \textit{external product} of $\mathbf{E}$ and $\mathbf{F}$ as the smooth vector bundle over $M\times N$ given by \[
\mathbf{E}\boxtimes \mathbf{F}:=\pi_M^*\mathbf{E}\otimes \pi_N^*\mathbf{F}.
\] In particular, for each $(x,y)\in M\times N$, its fiber is $(\mathbf{E}\boxtimes \mathbf{F})_{(x,y)}=\mathbf{E}_x\otimes \mathbf{F}_y$. \end{definition}

\begin{remark}\label{obs:distribuciones-en-haces-isomorfismo-canonico-efecto-fibra-fibra-ya} There is a canonical isomorphism \[
(\mathbf{E}\boxtimes \mathbf{F})^*\cong \mathbf{E}^*\boxtimes \mathbf{F}^*.
\] Fiberwise, we have $(\mathbf{E}_x\otimes \mathbf{F}_y)^*\cong \mathbf{E}_x^*\otimes \mathbf{F}_y^*$, since the vector spaces involved are finite dimensional. These fiberwise isomorphisms are compatible with local trivializations and therefore define a smooth bundle isomorphism. \end{remark}

The external product allows an expression of the form $\phi(x,y)\in \mathbf{E}_x\otimes \mathbf{F}_y^*$ to be interpreted intrinsically as a section of $\mathbf{E}\boxtimes \mathbf{F}^*$ over $M\times N$. Before applying distributions in one variable, we specify an elementary isomorphism that will be used repeatedly. \begin{remark}[Convention] If $V$ is a finite-dimensional $\mathbb K$-vector space and $N$ is a smooth manifold with or without boundary, we identify $V$ with the trivial vector bundle $N\times V\to N$. Consequently, if $\mathbf{F}\to N$ is a vector bundle, we write simply $V\otimes \mathbf{F}$ for the tensor product of the trivial bundle $N\times V$ with $\mathbf{F}$. With this convention, \[
(V\otimes \mathbf{F})_y\cong V\otimes \mathbf{F}_y,
\qquad y\in N.
\] \end{remark} \begin{proposition}\label{prop: isomorfismo tensor finito secciones} Let $N$ be a smooth manifold with or without boundary, let $\mathbf{F}\longrightarrow N$ be a smooth vector bundle, let $V$ be a finite-dimensional $\mathbb{K}$-vector space, and let $K\subseteq N$ be compact. Then, for every $k\in\mathbb{N}_{0}$, the linear map \[
\Theta_K:V\otimes_{\pi} \Gamma_{K}^{k}(\mathbf{F})
\longrightarrow
\Gamma_{K}^{k}(V\otimes \mathbf{F}),
\qquad
\Theta_K(v\otimes \mathbf{s})(y):=v\otimes \mathbf{s}(y),
\] extended by linearity, is a linear isomorphism.

Moreover, if we fix a norm on $V$, a Riemannian metric $\mathbf{g}$ on $N$, and a bundle metric and connection on $\mathbf{F}$, then $\Theta_K$ is a topological isomorphism from $V\otimes_{\pi}\Gamma_K^k(\mathbf{F})$, equipped with the projective tensor norm induced by $\|\cdot\|_V$ and $\|\cdot\|_{C^k(K,\mathbf{F})}$, onto $\Gamma_K^k(V\otimes \mathbf{F})$ with the norm $\|\cdot\|_{C^k(K,V\otimes \mathbf{F})}$. \end{proposition}

\begin{proof} First prove that $\Theta_K$ is well defined. If $v\in V$ and $\mathbf{s}\in\Gamma_K^k(\mathbf{F})$, then, for each $y\in N$, $v\otimes \mathbf{s}(y)\in (V\otimes \mathbf{F})_y$, so the map $y\longmapsto v\otimes \mathbf{s}(y)$ defines a section of class $C^k$ of the bundle $V\otimes \mathbf{F}$. Moreover, $\supp(\Theta_K(v\otimes \mathbf{s}))
\subseteq
\supp(\mathbf{s})
\subseteq
K$. Therefore $\Theta_K(v\otimes \mathbf{s})\in\Gamma_K^k(V\otimes \mathbf{F})$. Extending by linearity gives a linear map \[
\Theta_K\colon V\otimes\Gamma_K^k(\mathbf{F})\longrightarrow\Gamma_K^k(V\otimes \mathbf{F}).
\]

Let us check surjectivity. Take a basis $\{e_1,\dots,e_r\}$ of $V$, and let $\{\varepsilon^1,\dots,\varepsilon^r\}$ be the dual basis of $V^*$. If $\boldsymbol{\sigma}\in \Gamma_K^k(V\otimes \mathbf{F})$, define $\boldsymbol{\sigma}_i(y):=(\varepsilon^i\otimes \operatorname{id}_{\mathbf{F}_y})(\boldsymbol{\sigma}(y))$ for $i\in\{1,\dots,r\}$.

Since $\varepsilon^i\otimes \operatorname{id}_{\mathbf{F}}$ defines a smooth bundle homomorphism $V\otimes \mathbf{F}\longrightarrow \mathbf{F}$, it follows that $\boldsymbol{\sigma}_i$ is of class $C^k$. Moreover, if $y\notin K$, then $\boldsymbol{\sigma}(y)=0$ and hence $\boldsymbol{\sigma}_i(y)=0$. Thus $\boldsymbol{\sigma}_i\in \Gamma_K^k(\mathbf{F})$.

For each $y\in N$, using the basis $\{e_i\}$ gives $\boldsymbol{\sigma}(y)=\displaystyle\sum_{i=1}^r e_i\otimes \boldsymbol{\sigma}_i(y)$. Therefore $\boldsymbol{\sigma}=\Theta_K\left(\displaystyle\sum_{i=1}^r e_i\otimes \boldsymbol{\sigma}_i\right)$. This proves that $\Theta_K$ is surjective.

Now prove injectivity. Suppose that \[
\Theta_K\left(\displaystyle\sum_{i=1}^{r} e_i\otimes \mathbf{s}_i\right)=0,\qquad \mathbf{s}_i\in \Gamma_K^k(\mathbf{F}).
\] Then, for every $y\in N$, $\displaystyle\sum_{i=1}^{r}e_i\otimes \mathbf{s}_i(y)=0$ in $V\otimes \mathbf{F}_y$. Composing with $\varepsilon^j\otimes \operatorname{id}_{\mathbf{F}_y}$ gives $\mathbf{s}_j(y)=0$ for every $y\in N$. Since this holds for every $j\in\{1,\dots,r\}$, it follows that $\mathbf{s}_j=0$ for every $j$. Therefore $\displaystyle\sum_{i=1}^{r}e_i\otimes \mathbf{s}_i=0$ in $V\otimes \Gamma_K^k(\mathbf{F})$, and $\Theta_K$ is injective.

Finally, prove the topological assertion. The connection on $\mathbf{F}$ and the trivial connection on $N\times V$ induce a connection on $V\otimes \mathbf{F}$. If $\displaystyle \boldsymbol{\sigma}=\displaystyle\sum_{i=1}^{r}e_i\otimes \mathbf{s}_i$, then, for $0\leq j\leq k$, $\displaystyle \nabla^j\boldsymbol{\sigma}=\displaystyle\sum_{i=1}^{r}e_i\otimes\nabla^j \mathbf{s}_i$. Since $V$ is finite dimensional, there exist constants $A,B>0$, depending on the basis and fixed metrics but independent of the sections, such that for every $y\in K$ and $0\leq j\leq k$, \[
A\max_{1\leq i\leq r}|\nabla^j \mathbf{s}_i(y)|_{\mathbf{g},\mathbf{h}_{\mathbf{F}}}
\leq
\left|\displaystyle\sum_{i=1}^{r} e_i\otimes \nabla^j \mathbf{s}_i(y)\right|_{\mathbf{g},\mathbf{h}_{V\otimes \mathbf{F}}}
\leq
B\max_{1\leq i\leq r}|\nabla^j \mathbf{s}_i(y)|_{\mathbf{g},\mathbf{h}_{\mathbf{F}}}.
\] Taking the supremum over $y\in K$ and the maximum over $0\leq j\leq k$ gives constants $C_1,C_2>0$ such that \[
C_1\max_{1\leq i\leq r}\|\mathbf{s}_i\|_{C^k(K,\mathbf{F})}
\leq
\left\|\Theta_K\!\left(\displaystyle\sum_{i=1}^{r}e_i\otimes \mathbf{s}_i\right)
\right\|_{C^k(K,V\otimes \mathbf{F})}
\leq
C_2\max_{1\leq i\leq r}\|\mathbf{s}_i\|_{C^k(K,\mathbf{F})}.
\] For the fixed basis $(e_i)$, the quantity $\displaystyle \max_{1\leq i\leq r}\|\mathbf{s}_i\|_{C^k(K,\mathbf{F})}$ is equivalent to the projective norm of $\displaystyle \sum_{i=1}^{r} e_i\otimes \mathbf{s}_i$. These inequalities show that $\Theta_K$ and $\Theta_K^{-1}$ are continuous. \end{proof}

\begin{corollary}\label{cor: isomorfismo Dk tensor} Let $N$ be a smooth manifold with or without boundary, let $\mathbf{F}\longrightarrow N$ be a smooth vector bundle, and let $V$ be a finite-dimensional $\mathbb{K}$-vector space. Then, for every $k\in\mathbb{N}_0$, the linear map \[
\Theta:V\otimes \mathcal{D}^k(N,\mathbf{F})
\longrightarrow
\mathcal{D}^k(N,V\otimes \mathbf{F}),
\qquad
\Theta(v\otimes \mathbf{s})(y):=v\otimes \mathbf{s}(y),
\] extended by linearity, is a continuous linear isomorphism with continuous inverse. \end{corollary}

\begin{proof} Recall that $\mathcal{D}^k(N,\mathbf{F})=\Gamma_c^k(\mathbf{F})$ and that this space is equipped with the inductive topology \[
\mathcal{D}^k(N,\mathbf{F})=\varinjlim_{\ell\in\mathbb{N}}\Gamma_{K_\ell}^k(\mathbf{F}),
\] where $(K_\ell)_{\ell\in\mathbb{N}}$ is an increasing exhaustion of $N$ by compact sets. Likewise, \[
\mathcal{D}^k(N,V\otimes \mathbf{F})=\varinjlim_{\ell\in\mathbb{N}}\Gamma_{K_\ell}^k(V\otimes \mathbf{F}).
\]

For each $\ell$, Proposition~\ref{prop: isomorfismo tensor finito secciones} provides a continuous linear isomorphism with continuous inverse \[
\Theta_{K_\ell}:V\otimes \Gamma_{K_\ell}^k(\mathbf{F})
\longrightarrow
\Gamma_{K_\ell}^k(V\otimes \mathbf{F}).
\] Moreover, if $K_\ell\subseteq K_{\ell'}$, the corresponding maps are compatible with the natural inclusions, since all are defined by the same formula $v\otimes \mathbf{s}\mapsto (y\mapsto v\otimes \mathbf{s}(y))$.

Thus the isomorphisms $\Theta_{K_\ell}$ induce a linear map \[
\Theta:V\otimes \mathcal{D}^k(N,\mathbf{F})
\longrightarrow
\mathcal{D}^k(N,V\otimes \mathbf{F}).
\]

Let us prove that this map is bijective. If $\mathbf{s}\in\mathcal{D}^k(N,\mathbf{F})$, there exists $\ell$ such that $\mathbf{s}\in \Gamma_{K_\ell}^k(\mathbf{F})$; by the preceding proposition, $\Theta(v\otimes \mathbf{s})\in \Gamma_{K_\ell}^k(V\otimes \mathbf{F})\subseteq \mathcal{D}^k(N,V\otimes \mathbf{F})$. This shows that $\Theta$ is well defined.

Now let $\boldsymbol{\sigma}\in \mathcal{D}^k(N,V\otimes \mathbf{F})$. There exists $\ell$ such that $\boldsymbol{\sigma}\in \Gamma_{K_\ell}^k(V\otimes \mathbf{F})$. By Proposition~\ref{prop: isomorfismo tensor finito secciones}, there exist $\mathbf{s}_1,\dots,\mathbf{s}_r\in\Gamma_{K_\ell}^k(\mathbf{F})$ such that \[
\boldsymbol{\sigma}=\displaystyle\sum_{i=1}^{r}e_i\otimes \mathbf{s}_i.
\] Since $\mathbf{s}_i\in\mathcal{D}^k(N,\mathbf{F})$, we obtain \[
\boldsymbol{\sigma}=\Theta\left(\displaystyle\sum_{i=1}^{r}e_i\otimes \mathbf{s}_i\right).
\] This proves surjectivity. Injectivity follows in the same way by choosing a compact set containing the supports of all sections involved and applying injectivity of $\Theta_K$ on that compact set.

Finally, continuity of $\Theta$ and its inverse follows from the definition of the inductive topology, since it suffices to check continuity on each space $\Gamma_{K_\ell}^k$, where it was already proved in Proposition~\ref{prop: isomorfismo tensor finito secciones}. This proves that $\Theta$ is a topological isomorphism. \end{proof}

We apply this corollary with $V=\mathbf{E}_x$ and the bundle $\mathbf{F}^*$ over $N$. Thus, for each $x\in M$, there is a canonical isomorphism \[
\mathbf{E}_x\otimes \mathcal{D}^k(N,\mathbf{F}^*)\cong
\mathcal{D}^k(N,\mathbf{E}_x\otimes \mathbf{F}^*).
\] This justifies the action of a distribution in the variable $y$ on sections of $\mathbf{E}_x\otimes \mathbf{F}^*$.

\begin{lemma}\label{lem: distribucion dependiente de parametro} Let $(M,\mathbf{g})$ and $(N,\widetilde{\mathbf{g}})$ be Riemannian manifolds with or without boundary, equipped with their respective Riemann--Lebesgue measures $\lambda_{\mathbf{g}}$ and $\lambda_{\widetilde{\mathbf{g}}}$. Let $\mathbf{E}\longrightarrow M$ and $\mathbf{F}\longrightarrow N$ be smooth vector bundles, let $K\subseteq N$ be compact, and let $\boldsymbol{\phi}\in \Gamma^{k}(M\times N,\mathbf{E}\boxtimes \mathbf{F}^*)$ be a section such that $\supp(\boldsymbol{\phi})\subseteq M\times K$.

Let $m\in\mathbb N_0$, with $m<k$, and let $\mathbf T\in\mathcal D'(N,\mathbf F,\mathbb K)$. Suppose there are a compact set $L\subseteq N$, with $K\subseteq\operatorname{Int}_N L$, and a constant $C>0$ such that
\[
|\mathbf T(\boldsymbol\psi)|\leq C\|\boldsymbol\psi\|_{C^m(L,\mathbf F^*)},
\qquad \boldsymbol\psi\in\Gamma_L(\mathbf F^*).
\]
In particular, this hypothesis holds if $\mathbf T$ has order at most $m$. Then the map $\mathbf{f}\colon M\longrightarrow \mathbf{E}$ given by $\mathbf{f}(x):=\mathbf{T}(\boldsymbol{\phi}(x,\cdot))$ defines a section in $\Gamma^{k-m}(M,\mathbf{E})$ supported in the projection of $\supp(\boldsymbol{\phi})$ onto the first factor, that is, \[
\supp(\mathbf{f})\subseteq
\{x\in M\mid \exists y\in K \text{ such that } (x,y)\in\supp(\boldsymbol{\phi})\}.
\]

Moreover, if $P$ is a linear differential operator of order $\leq k-m$ acting on sections of $\mathbf{E}$, then \[
P\mathbf{f}=\mathbf{T}(P_x\boldsymbol{\phi}(x,\cdot)),
\] where $P_x$ indicates that $P$ acts in the variable $x$.

In particular, if $\boldsymbol{\phi}$ is smooth and compactly supported and $\mathbf{T}\in\mathcal{D}'(N,\mathbf{F},\mathbb{K})$ is an arbitrary distribution, then $\mathbf{f}$ is a compactly supported smooth section of $\mathbf{E}$. \end{lemma}

\begin{proof}
The order hypothesis allows action on sections of finite regularity with fixed support. Indeed, by Lemma~\ref{lem:densidad-soporte-aumentado-Cm}, each $\boldsymbol\psi\in\Gamma_K^m(\mathbf F^*)$ can be approximated in the $C^m$ norm by smooth sections supported in $L$. The stated bound makes their images under $\mathbf T$ converge and shows that the limit is independent of the approximation. This gives a continuous linear map on $\Gamma_K^m(\mathbf F^*)$ with the same bound. The construction is also independent of the auxiliary compact set: approximation can be carried out inside a common compact neighborhood of $K$.

For each $x\in M$, the isomorphism in Corollary~\ref{cor: isomorfismo Dk tensor} identifies $\boldsymbol\phi(x,\cdot)$ with an element of $\mathbf E_x\otimes\Gamma_K^m(\mathbf F^*)$. Define
\[
\mathbf f(x):=(\operatorname{id}_{\mathbf E_x}\otimes\mathbf T)
                   (\boldsymbol\phi(x,\cdot)).
\]
This definition does not depend on a basis of $\mathbf E_x$.

Work in a chart of $M$ and a frame of $\mathbf E$ of rank $r$. Put $\ell=k-m$ and, for $|\alpha|\leq\ell$, consider
\[
\Phi_\alpha(x):=D_x^\alpha\boldsymbol\phi(x,\cdot)
\in\Gamma_K^m(\mathbb K^r\otimes\mathbf F^*).
\]
These maps are continuous in the $C^m(K)$ norm. To see this, cover $K$ by finitely many trivializing charts. Every mixed derivative involved has total order at most $|\alpha|+m\leq k$; these derivatives are continuous and uniformly continuous on the product of a compact neighborhood of $x$ with $K$. Coordinate norms are equivalent to the $C^m(K)$ norm, so uniform continuity proves the assertion.

If $|\alpha|<\ell$ and the segment $x+th$, $0\leq t\leq1$, remains in the chart, the fundamental theorem of calculus, applied also to derivatives in $y$, gives
\[
\begin{split}
&\left\|\Phi_\alpha(x+h)-\Phi_\alpha(x)
          -\sum_{i=1}^nh_i\Phi_{\alpha+e_i}(x)\right\|_{C^m(K)}\\
&\qquad\leq\sum_{i=1}^n|h_i|\sup_{0\leq t\leq1}
\|\Phi_{\alpha+e_i}(x+th)-\Phi_{\alpha+e_i}(x)\|_{C^m(K)}.
\end{split}
\]
If $n=0$, the chart consists of one point and regularity is immediate. For $n\geq1$ and $h\ne0$, dividing by $|h|$ and using $\sum_{i=1}^n|h_i|\leq\sqrt n\,|h|$ gives
\[
\frac{\left\|\Phi_\alpha(x+h)-\Phi_\alpha(x)
          -\sum_{i=1}^nh_i\Phi_{\alpha+e_i}(x)\right\|_{C^m(K)}}{|h|}
\leq\sqrt n\max_{1\leq i\leq n}\sup_{0\leq t\leq1}
\|\Phi_{\alpha+e_i}(x+th)-\Phi_{\alpha+e_i}(x)\|_{C^m(K)}.
\]
Given $\varepsilon>0$, continuity of these $n$ maps allows us to choose $\delta>0$ so that each difference on the right is less than $\displaystyle\frac{\varepsilon}{\sqrt n}$ whenever $|h|<\delta$. Indeed, $|th|\leq|h|$ for every $t\in[0,1]$. The quotient is then less than $\varepsilon$ and therefore tends to zero. This proves norm differentiability of each $\Phi_\alpha$ of order below $\ell$, with derivative given by the continuous maps $\Phi_{\alpha+e_i}$. By induction, $x\mapsto\Phi_0(x)$ is $C^\ell$ with values in this Banach space. In boundary charts we use convex half-balls and segments contained in the half-space; the same bound gives inward derivatives and their continuity up to the boundary. The finite reflection extension described in the proof of Theorem~\ref{teo:diferenciabilidad-debil-fuerte-hilbert} applies to this Banach space as well and specifies $C^\ell$ regularity in those charts.

Continuity of $\operatorname{id}_{\mathbb K^r}\otimes\mathbf T$ allows this map to be applied to difference quotients. Thus
\[
D_x^\alpha\mathbf f(x)=\mathbf T(D_x^\alpha\boldsymbol\phi(x,\cdot)),
\qquad |\alpha|\leq\ell,
\]
and the derivatives are continuous. Hence $\mathbf f\in\Gamma^{k-m}(M,\mathbf E)$. The local expression of a differential operator $P$ as a sum of derivatives with smooth matrix coefficients, acting only on the factor $\mathbf E_x$, gives $P\mathbf f=\mathbf T(P_x\boldsymbol\phi)$.

The projection of $\supp(\boldsymbol\phi)$ onto $M$ is closed because the support is contained in $M\times K$ and $K$ is compact. Outside that projection there is a neighborhood $U_x$ on which $\boldsymbol\phi$ vanishes throughout $U_x\times N$. Consequently $\mathbf f$ vanishes on $U_x$, proving the support assertion.

If $\boldsymbol\phi$ is smooth and compactly supported and $\mathbf T$ is arbitrary, choose $K$ to be the projection of that support onto $N$ and a compact set $L$ whose interior contains $K$. Continuity of $\mathbf T$ on $\Gamma_L(\mathbf F^*)$ supplies a bound of finite order $m$. The argument applies for any $k>m$, so $\mathbf f$ is smooth. Its support is contained in the compact projection of $\supp(\boldsymbol\phi)$ onto $M$.
\end{proof}

The preceding result allows a distribution to be evaluated in one variable of a section defined on a product. The resulting section in the remaining variable loses at most as many derivatives as the order of the distribution.

\subsection{Tensor products of test spaces}

The topological description of test sections on $M\times N$ is obtained first with fixed supports, then by passage to the inductive limit. The first stage involves projective products of nuclear Fréchet spaces. The second involves the inductive product, which linearizes separately continuous bilinear maps on $LF$ spaces. Definitions and abstract results concerning these topologies are given in Section~\ref{ap:productos-tensoriales-topologicos-nucleares}.

\begin{lemma}[Localization in models with faces] \label{lem:nuclearidad-tensorizacion-modelos-con-frontera} Let $U\subseteq G_{m,r}:=\mathbb R^{m-r}\times[0,\infty)^r$ and $V\subseteq G_{n,s}$ be relatively open subsets. For every compact set $C\subseteq U$, the space $\mathcal D_C(U)$ of functions smooth up to the faces and supported in $C$ is nuclear Fréchet. If \[
C_0\subseteq\operatorname{Int}_{G_{m,r}}(C_1)\subseteq U,
\qquad
D_0\subseteq\operatorname{Int}_{G_{n,s}}(D_1)\subseteq V
\] are compact, there exists a continuous linear operator \[
\mathscr R:\mathcal D_{C_0\times D_0}(U\times V)
\longrightarrow
\mathcal D_{C_1}(U)\widehat\otimes_\pi\mathcal D_{D_1}(V)
\] whose composition with multiplication of functions is the identity. \end{lemma}

\begin{proof} For $f\in\mathcal D_C(U)$, its extension by zero to the complement of $U$ in $G_{m,r}$ is smooth: the support $C$ lies in the relatively open subset $U$. This operation crosses only the boundary of the open subset within the model. Next apply the operator $E_{m,r}$ from Lemma~\ref{lem:extension-suave-semiespacio-soporte-fijo}. We obtain a continuous map \[
J_C:\mathcal D_C(U)\longrightarrow
\mathcal D_{\widehat C}(\mathbb R^m),
\qquad J_Cf=E_{m,r}\widetilde f,
\] for a fixed compact set $\widehat C$. Restriction to $U$ is the inverse of $J_C$ on its image and is continuous for every derivative seminorm. This image is closed. Indeed, if $J_Cf_j\to h$ in $\mathcal D_{\widehat C}(\mathbb R^m)$, the restrictions converge with all derivatives to $f=h|_U$, which is supported in $C$; continuity of extension gives $J_Cf=h$. Thus $\mathcal D_C(U)$ is isomorphic to a closed subspace of a nuclear Fréchet space, and the first assertion follows from Corollary~\ref{cor:D-K-euclidiano-nuclear} and Proposition~\ref{prop:estabilidad-nuclear-subespacios-productos}.

To construct $\mathscr R$, first extend $\Phi\in\mathcal D_{C_0\times D_0}(U\times V)$ by zero to $G_{m,r}\times G_{n,s}$, then extend in each group of variables: \[
\widehat\Phi:=E_{m,r}^{x}E_{n,s}^{y}\widetilde\Phi.
\] The extension depends continuously on $\Phi$ in every seminorm and is supported in $A_0\times B_0$, for fixed compact sets $A_0\subseteq\mathbb R^m$ and $B_0\subseteq\mathbb R^n$. This follows by applying support control in each group while treating the other as a parameter. Choose compact sets $A_1,B_1$ containing $A_0,B_0$ in their respective interiors. Lemma~\ref{lem:tensorizacion-localizada-euclidiana} provides a continuous operator \[
R_{\mathrm{eucl}}:
\mathcal D_{A_0\times B_0}(\mathbb R^{m+n})
\longrightarrow
\mathcal D_{A_1}(\mathbb R^m)\widehat\otimes_\pi
\mathcal D_{B_1}(\mathbb R^n)
\] recovering each function upon application of the external product. Its construction by Fourier series on tori applies to $\widehat\Phi$, which is smooth across all faces.

Take $\chi\in\mathcal D_{C_1}(U)$ and $\eta\in\mathcal D_{D_1}(V)$ equal to one near $C_0$ and $D_0$. These functions are obtained by restricting Euclidean cutoffs: compactness allows each smaller compact set to be covered by finitely many balls whose intersections with the model lie in the relative interior of the larger compact set. The maps \[
P(f):=\chi(f|_U),\qquad Q(g):=\eta(g|_V)
\] are continuous into $\mathcal D_{C_1}(U)$ and $\mathcal D_{D_1}(V)$ by the Leibniz rule. Define \[
\mathscr R\Phi:=(P\widehat\otimes Q)
R_{\mathrm{eucl}}\widehat\Phi.
\] All operations are continuous and, by the restriction identity for the extensions, the product of their factors is \[
\chi(x)\eta(y)\widehat\Phi|_{U\times V}
=\chi(x)\eta(y)\Phi=\Phi.
\] This proves the second assertion. \end{proof}

\begin{proposition}[Nuclearity of sections with fixed support] \label{prop:Gamma-K-haz-nuclear} Let $M$ be a smooth manifold with corners, let $\mathbf{E}\to M$ be a smooth vector bundle of finite rank, and let $K\subseteq M$ be compact. Then $\Gamma_K(\mathbf{E})$ is a nuclear Fréchet space. \end{proposition}

\begin{proof} By Proposition~\ref{prop:prueba-esquinas-Frechet-LF}, only nuclearity remains to be proved. Choose coordinate and trivializing open subsets $U_1,\dots,U_s$ covering $K$, with $U_j$ relatively compact, and functions $\chi_j\in C_c^\infty(U_j)$ such that $\displaystyle\sum_{j=1}^s\chi_j=1$ on a neighborhood of $K$. Fix a chart $\kappa_j\colon U_j\longrightarrow\kappa_j(U_j)$ whose image is a relatively open subset of some $G_{\dim M,r_j}$, and a local frame $(\mathbf{e}_{j,1},\dots,\mathbf{e}_{j,r})$ of $\mathbf{E}$ on $U_j$, where $r=\operatorname{rank}(\mathbf{E})$.

For $\mathbf{u}\in\Gamma_K(\mathbf{E})$, write \[
(\chi_j\mathbf{u})|_{U_j}
=
\displaystyle\sum_{a=1}^r u_j^a \mathbf{e}_{j,a}
\] and set $K_j:=\kappa_j(K\cap\operatorname{supp}\chi_j)$. The map \[
I\colon \Gamma_K(\mathbf{E})\longrightarrow
\prod_{j=1}^s\mathcal D_{K_j}(\kappa_j(U_j))^r,
\qquad
I(\mathbf{u}):=\bigl(u_j^a\circ\kappa_j^{-1}\bigr)_{j,a},
\] is linear and continuous. Define a map in the opposite direction by \[
S\bigl((v_j^a)_{j,a}\bigr)
:=
\sum_{j=1}^s
\widetilde{
\left(
\sum_{a=1}^r(v_j^a\circ\kappa_j)\mathbf{e}_{j,a}
\right)}^{\,M}.
\] The tilde with superscript $M$ denotes extension by zero from $U_j$ to $M$. Each section in the sum is supported in the compact set $\kappa_j^{-1}(K_j)\Subset U_j$, so this extension is smooth. Moreover, $\kappa_j^{-1}(K_j)\subseteq K$, so $S$ takes values in $\Gamma_K(\mathbf{E})$. Local derivative expressions show that $S$ is continuous: each $C^k$ seminorm of a summand is bounded by a finite sum of $C^k$ seminorms of its component functions.

Since $\displaystyle\sum_{j=1}^s\chi_j=1$ on a neighborhood of $K$ and $\mathbf{u}$ vanishes outside $K$, we have $S(I(\mathbf{u}))=\mathbf{u}$ for every $\mathbf{u}\in\Gamma_K(\mathbf{E})$. Consequently, $I$ is a topological embedding, and $I\circ S$ is a continuous projection onto its image. This image is therefore a complemented subspace of the finite product.

Each space $\mathcal D_{K_j}(\kappa_j(U_j))$ is nuclear by Lemma~\ref{lem:nuclearidad-tensorizacion-modelos-con-frontera}. Proposition~\ref{prop:estabilidad-nuclear-subespacios-productos} shows that the finite product and its complemented subspace $I(\Gamma_K(\mathbf{E}))$ are nuclear. Since $I$ is a topological isomorphism onto its image, $\Gamma_K(\mathbf{E})$ is nuclear. \end{proof}

\begin{theorem}[Localized tensor product representation of sections] \label{teo:tensorizacion-localizada-secciones} Let $M$ and $N$ be smooth manifolds with corners, and let $\mathbf{E}\to M$ and $\mathbf{F}\to N$ be smooth vector bundles of finite rank. Let the following sets be compact: \[
K_0\subseteq\operatorname{Int}(K_1)\subseteq M,
\qquad
L_0\subseteq\operatorname{Int}(L_1)\subseteq N.
\] The external product induces a topological embedding \[
T_{K_1,L_1}:
\Gamma_{K_1}(\mathbf{E})\widehat\otimes_\pi\Gamma_{L_1}(\mathbf{F})
\longrightarrow
\Gamma_{K_1\times L_1}(\mathbf{E}\boxtimes \mathbf{F}).
\] Its image contains $\Gamma_{K_0\times L_0}(\mathbf{E}\boxtimes \mathbf{F})$, and the restriction of $T_{K_1,L_1}^{-1}$ to this subspace is continuous. \end{theorem}

\begin{proof} Local derivative estimates for $u\boxtimes v$ show that the bilinear map \[
\Gamma_{K_1}(\mathbf{E})\times\Gamma_{L_1}(\mathbf{F})
\longrightarrow
\Gamma_{K_1\times L_1}(\mathbf{E}\boxtimes \mathbf{F})
\] is continuous. The universal property of the projective product and completeness of the target give the continuous operator $T_{K_1,L_1}$.

We first construct a continuous preimage for sections supported in $K_0\times L_0$. Choose coordinate and trivializing open subsets \[
U_1,\dots,U_A\Subset\operatorname{Int}(K_1),
\qquad
V_1,\dots,V_B\Subset\operatorname{Int}(L_1),
\] covering $K_0$ and $L_0$. Let $(\mathbf{e}_{i,1},\dots,\mathbf{e}_{i,r})$ and $(\mathbf{f}_{j,1},\dots,\mathbf{f}_{j,s})$ be the local frames. Take cutoff functions $\chi_i\in C_c^\infty(U_i)$ and $\eta_j\in C_c^\infty(V_j)$ such that \[
\sum_{i=1}^A\chi_i=1
\quad\text{on a neighborhood of }K_0,
\qquad
\sum_{j=1}^B\eta_j=1
\quad\text{on a neighborhood of }L_0.
\] Denote the charts by $\kappa_i\colon U_i\longrightarrow\kappa_i(U_i)$ and $\lambda_j\colon V_j\longrightarrow\lambda_j(V_j)$. Choose compact sets $\widetilde C_i\Subset\kappa_i(U_i)$ and $\widetilde D_j\Subset\lambda_j(V_j)$ whose interiors contain, respectively, \[
C_i:=\kappa_i(K_0\cap\operatorname{supp}(\chi_i)),
\qquad
D_j:=\lambda_j(L_0\cap\operatorname{supp}(\eta_j)).
\]

Let $\boldsymbol{\Phi}\in\Gamma_{K_0\times L_0}(\mathbf{E}\boxtimes \mathbf{F})$. On $U_i\times V_j$, write \[
(\chi_i\boxtimes\eta_j)\boldsymbol{\Phi}
=
\sum_{a=1}^r\sum_{b=1}^s
\bigl(\Phi_{ij}^{ab}\circ(\kappa_i\times\lambda_j)\bigr)
 \mathbf{e}_{i,a}\boxtimes \mathbf{f}_{j,b},
\] where \[
\Phi_{ij}^{ab}\in
\mathcal D_{C_i\times D_j}
\bigl(\kappa_i(U_i)\times\lambda_j(V_j)\bigr).
\] Extraction of each coefficient is continuous in the fixed-support topologies. Lemma~\ref{lem:nuclearidad-tensorizacion-modelos-con-frontera} gives a continuous operator \[
\mathscr R_{ij}:
\mathcal D_{C_i\times D_j}
\bigl(\kappa_i(U_i)\times\lambda_j(V_j)\bigr)
\longrightarrow
\mathcal D_{\widetilde C_i}(\kappa_i(U_i))
\widehat\otimes_\pi
\mathcal D_{\widetilde D_j}(\lambda_j(V_j))
\] whose composition with multiplication of functions is the identity.

For $h\in\mathcal D_{\widetilde C_i}(\kappa_i(U_i))$, the section $(h\circ\kappa_i)\mathbf{e}_{i,a}$ is supported in $\kappa_i^{-1}(\widetilde C_i)\Subset U_i$. Denote by \[
\widetilde{(h\circ\kappa_i)\mathbf{e}_{i,a}}^{\,M}
\] its extension by zero from $U_i$ to $M$. If $k\in\mathcal D_{\widetilde D_j}(\lambda_j(V_j))$, write \[
\widetilde{(k\circ\lambda_j)\mathbf{f}_{j,b}}^{\,N}
\] for the extension by zero from $V_j$ to $N$. Define \[
P_{i,a}(h):=\widetilde{(h\circ\kappa_i)\mathbf{e}_{i,a}}^{\,M},
\qquad
Q_{j,b}(k):=\widetilde{(k\circ\lambda_j)\mathbf{f}_{j,b}}^{\,N}.
\] These extensions are smooth, and the maps $P_{i,a}$ and $Q_{j,b}$ are continuous. Moreover, they take values in $\Gamma_{K_1}(\mathbf{E})$ and $\Gamma_{L_1}(\mathbf{F})$ because their supports lie in $\operatorname{Int}(K_1)$ and $\operatorname{Int}(L_1)$. By the universal property of the projective product, $P_{i,a}\widehat\otimes Q_{j,b}$ is continuous. Set \[
R_{K_0,L_0}^{K_1,L_1}\boldsymbol{\Phi}
:=
\sum_{i=1}^A\sum_{j=1}^B
\sum_{a=1}^r\sum_{b=1}^s
(P_{i,a}\widehat\otimes Q_{j,b})
\mathscr R_{ij}(\Phi_{ij}^{ab}).
\] The sum is finite, so $R_{K_0,L_0}^{K_1,L_1}$ is continuous. Applying $T_{K_1,L_1}$, each term recovers the corresponding coefficient of $(\chi_i\boxtimes\eta_j)\boldsymbol{\Phi}$. The partition-of-unity identities give \[
T_{K_1,L_1}R_{K_0,L_0}^{K_1,L_1}\boldsymbol{\Phi}
=
\sum_{i=1}^A\sum_{j=1}^B
(\chi_i\boxtimes\eta_j)\boldsymbol{\Phi}
=
\boldsymbol{\Phi}.
\]

It remains to prove that $T_{K_1,L_1}$ is a topological embedding. Cover $K_1$ and $L_1$ by finitely many trivializing charts $(U_\alpha,\kappa_\alpha)$ and $(V_\beta,\lambda_\beta)$. Choose cutoff functions $\zeta_\alpha\in C_c^\infty(U_\alpha)$ and $\theta_\beta\in C_c^\infty(V_\beta)$ whose sums equal one on neighborhoods of $K_1$ and $L_1$. Taking the components of $\zeta_\alpha u$ and $\theta_\beta v$ in the local frames, first extend them by zero within their coordinate models and then, by Lemma~\ref{lem:extension-suave-semiespacio-soporte-fijo}, to the full Euclidean space. In interior charts, the latter operator is the identity. We obtain topological embeddings \[
I_{\mathbf{E}}:\Gamma_{K_1}(\mathbf{E})
\hookrightarrow
\prod_{c=1}^{C}\mathcal D_{A_c}(\mathbb R^{\dim M}),
\qquad
I_{\mathbf{F}}:\Gamma_{L_1}(\mathbf{F})
\hookrightarrow
\prod_{d=1}^{D}\mathcal D_{B_d}(\mathbb R^{\dim N}).
\] The index $c$ ranges over pairs consisting of a chart of the first factor and a frame component, and $d$ ranges over the analogous pairs for the second factor. Let the compact supports before extension be \[
A_c^0=\kappa_\alpha(K_1\cap\operatorname{supp}(\zeta_\alpha)),
\qquad
B_d^0=\lambda_\beta(L_1\cap\operatorname{supp}(\theta_\beta)),
\] repeated once for each component. Choose $A_c$ and $B_d$ containing the fixed supports of their smooth extensions. For example, if $u_\alpha^a$ is the component of $\zeta_\alpha u$ in $U_\alpha$, the corresponding coordinate of $I_{\mathbf{E}}u$ is $E_\alpha\bigl(\widetilde{u_\alpha^a\circ\kappa_\alpha^{-1}}\bigr)$, where the tilde denotes extension by zero from $\kappa_\alpha(U_\alpha)$ to its model $G_{\dim M,r_\alpha}$, and $E_\alpha=E_{\dim M,r_\alpha}$, with $r_\alpha\in\{0,\ldots,\dim M\}$. The extension with a tilde is smooth because $\operatorname{supp}(\zeta_\alpha)\Subset U_\alpha$, and $E_\alpha$ smoothly extends its jets across the geometric boundary. The same convention, with operators $E_\beta=E_{\dim N,r_\beta}$, where $r_\beta\in\{0,\ldots,\dim N\}$, is used for $I_{\mathbf{F}}$.

That $I_{\mathbf{E}}$ and $I_{\mathbf{F}}$ are embeddings is checked with the $C^k$ seminorms: components are bounded by seminorms of the section, while the cutoff identity bounds the section seminorms by a finite sum of seminorms of its components. The same construction, using $\zeta_\alpha\boxtimes\theta_\beta$ and the extension $E_\alpha^xE_\beta^y$ of each localized component, gives an embedding \[
I_{\mathbf{E}\boxtimes \mathbf{F}}:
\Gamma_{K_1\times L_1}(\mathbf{E}\boxtimes \mathbf{F})
\hookrightarrow
\prod_{c,d}
\mathcal D_{A_c\times B_d}(\mathbb R^{\dim M+\dim N}).
\] Support control in each group of variables allows us to use precisely the compact sets $A_c\times B_d$. The product identity in the extension lemma ensures, for elementary sections, that extending a component of $u\boxtimes v$ is the same as multiplying the extended components of $u$ and $v$. By continuity, the same identity holds in the completed tensor products.

By Proposition~\ref{prop:Gamma-K-haz-nuclear}, the two domain spaces are nuclear. The finite products in the codomain are nuclear as well. Theorem~\ref{teo:nuclearidad-pi-igual-epsilon} identifies the topologies $\pi$ and $\varepsilon$, and Proposition~\ref{prop:caracter-inyectivo-epsilon} shows that $I_{\mathbf{E}}\widehat\otimes I_{\mathbf{F}}$ is a topological embedding. After applying Proposition~\ref{prop:tensor-proyectivo-productos-finitos}, the external product becomes, in each coordinate $(c,d)$, the scalar map \[
\mathcal D_{A_c}(\mathbb R^{\dim M})\widehat\otimes_\pi
\mathcal D_{B_d}(\mathbb R^{\dim N})
\longrightarrow
\mathcal D_{A_c\times B_d}(\mathbb R^{\dim M+\dim N}),
\] which is an embedding by Lemma~\ref{lem:tensorizacion-localizada-euclidiana}.

Let $(z_\gamma)$ be a net in $\Gamma_{K_1}(\mathbf{E})\widehat\otimes_\pi\Gamma_{L_1}(\mathbf{F})$, and suppose that $T_{K_1,L_1}z_\gamma\to0$. Commutativity of the localization maps implies that the image of $(I_{\mathbf{E}}\widehat\otimes I_{\mathbf{F}})z_\gamma$ under all the preceding scalar embeddings converges to zero. Since their finite product is an embedding, $(I_{\mathbf{E}}\widehat\otimes I_{\mathbf{F}})z_\gamma\to0$. The map $I_{\mathbf{E}}\widehat\otimes I_{\mathbf{F}}$ is also an embedding, so $z_\gamma\to0$. Applying the same argument to a constant net proves injectivity of $T_{K_1,L_1}$. Thus its inverse on the image is continuous. By injectivity, the operator constructed at the outset is the restriction of this inverse to $\Gamma_{K_0\times L_0}(\mathbf{E}\boxtimes \mathbf{F})$. \end{proof}

\begin{theorem}[Tensor product representation of test spaces]\label{teo: tensorizacion espacios prueba}\index{tensor product representation of test spaces@tensor product representation of test spaces} Let $M$ and $N$ be smooth manifolds with corners, and let $\mathbf{E}\to M$ and $\mathbf{F}\to N$ be smooth vector bundles of finite rank. The external product induces a canonical topological isomorphism \[
\mathcal D(M,\mathbf{E})\widehat\otimes_\iota\mathcal D(N,\mathbf{F})
\cong
\mathcal D(M\times N,\mathbf{E}\boxtimes \mathbf{F}).
\] In particular, finite linear combinations of sections $\mathbf{u}\boxtimes \mathbf{v}$ are dense in $\mathcal D(M\times N,\mathbf{E}\boxtimes \mathbf{F})$. \end{theorem}

\begin{proof} Choose exhaustions by compact sets $(K_j)_{j\in\mathbb N}$ and $(L_j)_{j\in\mathbb N}$ such that $K_j\subseteq\operatorname{Int}(K_{j+1})$ and $L_j\subseteq\operatorname{Int}(L_{j+1})$. Then $(K_j\times L_j)_{j\in\mathbb N}$ is an exhaustion of $M\times N$ with the same property.

First identify the topology on the algebraic tensor product. The $\iota$ topology on $\mathcal D(M,\mathbf{E})\otimes\mathcal D(N,\mathbf{F})$ agrees with the final locally convex topology for the maps \[
\Gamma_{K_p}(\mathbf{E})\otimes_\pi\Gamma_{L_q}(\mathbf{F})
\longrightarrow
\mathcal D(M,\mathbf{E})\otimes\mathcal D(N,\mathbf{F}),
\qquad p,q\in\mathbb N.
\] To check this, let $b\colon \mathcal D(M,\mathbf{E})\times\mathcal D(N,\mathbf{F})\longrightarrow Z$ be a bilinear map. If $b$ is separately continuous, its restriction to $\Gamma_{K_p}(\mathbf{E})\times\Gamma_{L_q}(\mathbf{F})$ is separately continuous as well. These two spaces are Fréchet, so Proposition~\ref{prop:bilineal-separadamente-continua-Frechet} makes the restriction continuous. Conversely, suppose all restrictions are continuous. If $\mathbf{u}\in\mathcal D(M,\mathbf{E})$, there exists $p$ such that $\mathbf{u}\in\Gamma_{K_p}(\mathbf{E})$. For each $q$, the map $\mathbf{v}\mapsto b(\mathbf{u},\mathbf{v})$ is continuous on $\Gamma_{L_q}(\mathbf{F})$. The final property of the $LF$ topology on $\mathcal D(N,\mathbf{F})$ implies continuity on the entire space. The same argument with the factors interchanged proves continuity in the first variable. The universal property of the inductive product, Proposition~\ref{prop:propiedad-universal-producto-inductivo}, gives equality of the topologies.

The external product is separately continuous. Indeed, if $\mathbf{v}\in\Gamma_{L_q}(\mathbf{F})$ is fixed, for each $p$ the local estimates $C^k$ show that $\mathbf{u}\longmapsto \mathbf{u}\boxtimes \mathbf{v}$ is continuous from $\Gamma_{K_p}(\mathbf{E})$ to $\Gamma_{K_p\times L_q}(\mathbf{E}\boxtimes \mathbf{F})$. The final property of the $LF$ topology gives continuity on $\mathcal D(M,\mathbf{E})$. The other factor is treated in the same way. Thus the external product linearizes to a continuous operator on $\mathcal D(M,\mathbf{E})\otimes_\iota\mathcal D(N,\mathbf{F})$. Since $\mathcal D(M\times N,\mathbf{E}\boxtimes \mathbf{F})$ is complete by Theorem~\ref{teo: propiedades LF}, this operator extends uniquely to \[
T:\mathcal D(M,\mathbf{E})\widehat\otimes_\iota\mathcal D(N,\mathbf{F})
\longrightarrow
\mathcal D(M\times N,\mathbf{E}\boxtimes \mathbf{F}).
\]

We construct the inverse of $T$. For $\boldsymbol{\Phi}\in\Gamma_{K_j\times L_j}(\mathbf{E}\boxtimes \mathbf{F})$, Theorem~\ref{teo:tensorizacion-localizada-secciones} gives a unique element \[
R_j\boldsymbol{\Phi}\in
\Gamma_{K_{j+1}}(\mathbf{E})\widehat\otimes_\pi\Gamma_{L_{j+1}}(\mathbf{F})
\] whose image under the external product is $\boldsymbol{\Phi}$. The map $R_j$ is continuous. Inclusion of the two fixed-support stages into the test spaces induces a continuous operator from this completed projective product to $\mathcal D(M,\mathbf{E})\widehat\otimes_\iota\mathcal D(N,\mathbf{F})$.

We must check independence of the stage. Suppose that $j<q$ and $\boldsymbol{\Phi}\in\Gamma_{K_j\times L_j}(\mathbf{E}\boxtimes \mathbf{F})$. Include $R_j\boldsymbol{\Phi}$ and $R_q\boldsymbol{\Phi}$ in \[
\Gamma_{K_{q+1}}(\mathbf{E})\widehat\otimes_\pi\Gamma_{L_{q+1}}(\mathbf{F}).
\] The inclusions used here are topological embeddings. Indeed, fixed-support spaces are nuclear by Proposition~\ref{prop:Gamma-K-haz-nuclear}, so $\pi=\varepsilon$, and the assertion follows from Proposition~\ref{prop:caracter-inyectivo-epsilon}. Both preimages have external-product image $\boldsymbol{\Phi}$. Injectivity in Theorem~\ref{teo:tensorizacion-localizada-secciones}, applied at stage $K_{q+1}\times L_{q+1}$, shows that they agree.

We may therefore define \[
R:\mathcal D(M\times N,\mathbf{E}\boxtimes \mathbf{F})
\longrightarrow
\mathcal D(M,\mathbf{E})\widehat\otimes_\iota\mathcal D(N,\mathbf{F}),
\qquad
R\boldsymbol{\Phi}:=R_j\boldsymbol{\Phi}
\] when $\operatorname{supp}\boldsymbol{\Phi}\subseteq K_j\times L_j$. Each restriction of $R$ to $\Gamma_{K_j\times L_j}(\mathbf{E}\boxtimes \mathbf{F})$ is continuous. The $LF$ topology on the domain is final with respect to these stages, so $R$ is continuous. By construction, $T\circ R$ is the identity.

Now let $\mathbf{u}\in\mathcal D(M,\mathbf{E})$ and $\mathbf{v}\in\mathcal D(N,\mathbf{F})$. Choose $j$ so that $\operatorname{supp}\mathbf{u}\subseteq K_j$ and $\operatorname{supp}\mathbf{v}\subseteq L_j$. The elements $R(\mathbf{u}\boxtimes \mathbf{v})$ and $\mathbf{u}\otimes \mathbf{v}$, viewed at the corresponding stage, have the same image under the external product. Localized injectivity gives \[
R(\mathbf{u}\boxtimes \mathbf{v})=\mathbf{u}\otimes \mathbf{v}.
\] Thus $R\circ T$ is the identity on the algebraic tensor product. This subspace is dense in its completion, and continuity of $R\circ T$ extends the identity to all of $\mathcal D(M,\mathbf{E})\widehat\otimes_\iota\mathcal D(N,\mathbf{F})$. Therefore $T$ and $R$ are inverse topological isomorphisms. The image of the algebraic tensor product under $T$ is the space of finite linear combinations of sections $\mathbf{u}\boxtimes \mathbf{v}$, so this space is dense in the test space of the product. \end{proof}

\begin{corollary}[Tensor decomposition with support control] \label{cor:descomposicion-tensorial-control-soporte} Let $M$ and $N$ be smooth manifolds with corners, let $\mathbf{E}\to M$ and $\mathbf{F}\to N$ be smooth vector bundles of finite rank, and let $K_0\subseteq\operatorname{Int}(K_1)\subseteq M$ and $L_0\subseteq\operatorname{Int}(L_1)\subseteq N$ be compact. Every section \[
\boldsymbol{\Phi}\in\Gamma_{K_0\times L_0}(\mathbf{E}^*\boxtimes \mathbf{F}^*)
\] has a unique preimage \[
u_{\boldsymbol{\Phi}}\in
\Gamma_{K_1}(\mathbf{E}^*)\widehat\otimes_\pi\Gamma_{L_1}(\mathbf{F}^*)
\] under the external-product map. Moreover, $\boldsymbol{\Phi}\mapsto u_{\boldsymbol{\Phi}}$ is continuous. \end{corollary}

\begin{proof} Apply Theorem~\ref{teo:tensorizacion-localizada-secciones} to the bundles $\mathbf{E}^*$ and $\mathbf{F}^*$. The external-product map at stage $K_1\times L_1$ is injective and its image contains $\Gamma_{K_0\times L_0}(\mathbf{E}^*\boxtimes \mathbf{F}^*)$. Thus each $\boldsymbol{\Phi}$ has a unique preimage $u_{\boldsymbol{\Phi}}$. Continuity of $\boldsymbol{\Phi}\mapsto u_{\boldsymbol{\Phi}}$ is precisely continuity of the inverse restricted to this subspace. \end{proof}

\subsection{The Schwartz kernel theorem}

A continuous operator can be studied as a distribution on the product of its source and target manifolds. In the presence of bundles, the kernel takes values in the appropriate tensor product; the following formulation makes this identification explicit.

\begin{theorem}[Schwartz kernel theorem for vector bundles]\label{teo: nucleo Schwartz haces}\index{Schwartz kernel theorem for vector bundles@Schwartz kernel theorem for vector bundles} Let $M$ and $N$ be smooth manifolds with corners, let $\mathbf{E}\to M$ and $\mathbf{F}\to N$ be smooth vector bundles, and let $W$ be a finite-dimensional $\mathbb{K}$-vector space.

Let \[
A\colon \mathcal{D}(N,\mathbf{F}^*)\longrightarrow \mathcal{D}'(M,\mathbf{E},W)
\] be a linear map, where the target is equipped with the weak-$*$ topology, or equivalently the pointwise convergence topology of Definition~\ref{def:distribuciones-en-haces-distribuciones}. The following conditions are equivalent: \begin{enumerate}[label=(\alph*)] \item $A$ is continuous; \item the bilinear map \[
B_A\colon \mathcal{D}(M,\mathbf{E}^*)\times \mathcal{D}(N,\mathbf{F}^*)\longrightarrow W,
\qquad
B_A(\boldsymbol{\phi},\boldsymbol{\psi}):=(A\boldsymbol{\psi})(\boldsymbol{\phi}),
\] is separately continuous; \item for any compact sets $K\subseteq M$ and $L\subseteq N$, the restriction of $B_A$ to $\Gamma_K(\mathbf{E}^*)\times\Gamma_L(\mathbf{F}^*)$ is continuous. \end{enumerate} When these conditions hold, there exists a unique distribution \[
\mathbf{K}_A\in \mathcal{D}'(M\times N,\mathbf{E}\boxtimes \mathbf{F},W)
\] such that \[
\mathbf{K}_A(\boldsymbol{\phi}\boxtimes \boldsymbol{\psi})=(A\boldsymbol{\psi})(\boldsymbol{\phi})
\] for every $\boldsymbol{\phi}\in \mathcal{D}(M,\mathbf{E}^*)$ and every $\boldsymbol{\psi}\in \mathcal{D}(N,\mathbf{F}^*)$.

Conversely, every distribution $\mathbf{K}\in \mathcal{D}'(M\times N,\mathbf{E}\boxtimes \mathbf{F},W)$ defines a linear map \[
A_{\mathbf{K}}\colon \mathcal{D}(N,\mathbf{F}^*)\longrightarrow \mathcal{D}'(M,\mathbf{E},W)
\] by \[
(A_{\mathbf{K}}\boldsymbol{\psi})(\boldsymbol{\phi}):=\mathbf{K}(\boldsymbol{\phi}\boxtimes \boldsymbol{\psi}),
\] and this map is continuous in the weak-$*$ topology. Moreover, the correspondences $A\mapsto \mathbf{K}_A$ and $\mathbf{K}\mapsto A_{\mathbf{K}}$ are linear and mutually inverse. In other words, there is a canonical linear bijection \[
\mathcal D'(M\times N,\mathbf{E}\boxtimes \mathbf{F},W)
\cong
\mathcal L\bigl(\mathcal D(N,\mathbf{F}^*),\mathcal D'(M,\mathbf{E},W)_{\sigma}\bigr),
\] where the subscript $\sigma$ indicates the weak-$*$ topology. \end{theorem}

\begin{proof} If $A$ is continuous in the topology of pointwise convergence, then, for each $\boldsymbol{\phi}\in\mathcal D(M,\mathbf{E}^*)$, the map \[
\boldsymbol{\psi}\longmapsto(A\boldsymbol{\psi})(\boldsymbol{\phi})
=
\operatorname{ev}_{\boldsymbol{\phi}}(A\boldsymbol{\psi})
\] is continuous. For each $\boldsymbol{\psi}$, the map $\boldsymbol{\phi}\mapsto(A\boldsymbol{\psi})(\boldsymbol{\phi})$ is continuous because $A\boldsymbol{\psi}$ is a distribution. This proves $(a)\Rightarrow(b)$.

Suppose $(b)$. The restriction of $B_A$ to $\Gamma_K(\mathbf{E}^*)\times\Gamma_L(\mathbf{F}^*)$ is separately continuous. Both factors are Fréchet spaces by Proposition~\ref{prop:prueba-esquinas-Frechet-LF}, and Proposition~\ref{prop:bilineal-separadamente-continua-Frechet} shows that the restriction is continuous. Therefore $(b)\Rightarrow(c)$.

Finally, suppose $(c)$ and fix $\boldsymbol{\phi}\in\mathcal D(M,\mathbf{E}^*)$. There exists a compact set $K\subseteq M$ such that $\boldsymbol{\phi}\in\Gamma_K(\mathbf{E}^*)$. For every compact set $L\subseteq N$, condition $(c)$ implies that $\boldsymbol{\psi}\mapsto B_A(\boldsymbol{\phi},\boldsymbol{\psi})$ is continuous on $\Gamma_L(\mathbf{F}^*)$. Since the topology on $\mathcal D(N,\mathbf{F}^*)$ is final with respect to these fixed-support spaces, the same map is continuous on $\mathcal D(N,\mathbf{F}^*)$. Thus $\operatorname{ev}_{\boldsymbol{\phi}}\circ A$ is continuous for every $\boldsymbol{\phi}$. The weak-$*$ topology on $\mathcal D'(M,\mathbf{E},W)$ is the initial topology for evaluations, so Proposition~\ref{prop: caracterizacion topologia inicial} proves continuity of $A$. This gives $(c)\Rightarrow(a)$.

Now suppose these conditions hold. By $(b)$ and the universal property of the inductive product, Proposition~\ref{prop:propiedad-universal-producto-inductivo}, there exists a unique continuous operator \[
\widetilde B_A:
\mathcal D(M,\mathbf{E}^*)\otimes_\iota\mathcal D(N,\mathbf{F}^*)
\longrightarrow W
\] such that $\widetilde B_A(\boldsymbol{\phi}\otimes\boldsymbol{\psi})=B_A(\boldsymbol{\phi},\boldsymbol{\psi})$. The space $W$ is complete because it is finite dimensional. The universal property of completion extends $\widetilde B_A$ uniquely to a continuous operator \[
\widehat B_A:
\mathcal D(M,\mathbf{E}^*)\widehat\otimes_\iota\mathcal D(N,\mathbf{F}^*)
\longrightarrow W.
\] Let \[
T:
\mathcal D(M,\mathbf{E}^*)\widehat\otimes_\iota\mathcal D(N,\mathbf{F}^*)
\longrightarrow
\mathcal D(M\times N,\mathbf{E}^*\boxtimes \mathbf{F}^*)
\] be the isomorphism in Theorem~\ref{teo: tensorizacion espacios prueba}. Define \[
\mathbf{K}_A:=\widehat B_A\circ T^{-1}.
\] This is a continuous linear operator on $\mathcal D(M\times N,\mathbf{E}^*\boxtimes \mathbf{F}^*)$. Through the isomorphism $(\mathbf{E}\boxtimes \mathbf{F})^*\cong \mathbf{E}^*\boxtimes \mathbf{F}^*$ in Remark~\ref{obs:distribuciones-en-haces-isomorfismo-canonico-efecto-fibra-fibra-ya}, it is identified with an element of $\mathcal D'(M\times N,\mathbf{E}\boxtimes \mathbf{F},W)$. Moreover, \[
\mathbf{K}_A(\boldsymbol{\phi}\boxtimes\boldsymbol{\psi})
=
\widehat B_A(\boldsymbol{\phi}\otimes\boldsymbol{\psi})
=
B_A(\boldsymbol{\phi},\boldsymbol{\psi})
=
(A\boldsymbol{\psi})(\boldsymbol{\phi}).
\]

For uniqueness, let $\mathbf{S}$ be another distribution satisfying the same identity. Then $\mathbf{S}-\mathbf{K}_A$ vanishes on every section $\boldsymbol{\phi}\boxtimes\boldsymbol{\psi}$ and, by linearity, on all finite linear combinations of them. These combinations are dense by Theorem~\ref{teo: tensorizacion espacios prueba}. Since $\mathbf{S}-\mathbf{K}_A$ is continuous, it vanishes on all of $\mathcal D(M\times N,\mathbf{E}^*\boxtimes \mathbf{F}^*)$. Therefore $\mathbf{S}=\mathbf{K}_A$.

Conversely, let $\mathbf{K}\in\mathcal D'(M\times N,\mathbf{E}\boxtimes \mathbf{F},W)$. For $\boldsymbol{\psi}\in\mathcal D(N,\mathbf{F}^*)$, define \[
(A_{\mathbf{K}}\boldsymbol{\psi})(\boldsymbol{\phi}):=\mathbf{K}(\boldsymbol{\phi}\boxtimes\boldsymbol{\psi}).
\] For fixed $\boldsymbol{\psi}$, the map $\boldsymbol{\phi}\mapsto\boldsymbol{\phi}\boxtimes\boldsymbol{\psi}$ is continuous, as proved in constructing the isomorphism of Theorem~\ref{teo: tensorizacion espacios prueba}. Its composition with $\mathbf{K}$ is therefore a continuous operator from $\mathcal D(M,\mathbf{E}^*)$ to $W$. This proves that $A_{\mathbf{K}}\boldsymbol{\psi}\in\mathcal D'(M,\mathbf{E},W)$. Linearity of $A_{\mathbf{K}}$ follows from linearity of $\mathbf{K}$.

To prove continuity of $A_{\mathbf{K}}$ into the weak-$*$ dual, fix $\boldsymbol{\phi}\in\mathcal D(M,\mathbf{E}^*)$. The map \[
\operatorname{ev}_{\boldsymbol{\phi}}\circ A_{\mathbf{K}}:
\boldsymbol{\psi}\longmapsto \mathbf{K}(\boldsymbol{\phi}\boxtimes\boldsymbol{\psi})
\] is continuous because $\boldsymbol{\psi}\mapsto\boldsymbol{\phi}\boxtimes\boldsymbol{\psi}$ is continuous. Every evaluation is continuous. By the universal property of the initial topology, $A_{\mathbf{K}}$ is continuous in the weak-$*$ topology.

Finally, \[
(A_{\mathbf{K}_A}\boldsymbol{\psi})(\boldsymbol{\phi})
=
\mathbf{K}_A(\boldsymbol{\phi}\boxtimes\boldsymbol{\psi})
=
(A\boldsymbol{\psi})(\boldsymbol{\phi}),
\] so $A_{\mathbf{K}_A}=A$. Starting with $\mathbf{K}$, the distributions $\mathbf{K}_{A_{\mathbf{K}}}$ and $\mathbf{K}$ agree on every elementary product. Density proved in Theorem~\ref{teo: tensorizacion espacios prueba} and continuity of both distributions imply that they agree on the entire test space. Therefore $\mathbf{K}_{A_{\mathbf{K}}}=\mathbf{K}$. The definitions also show that the two correspondences are linear. \end{proof}

\begin{example}[A kernel supported at a corner of the product]
\label{ej:nucleo-schwartz-delta-esquina}
Let $M=N=I=[0,1]$ and take trivial real scalar bundles.
Evaluation $\delta_0(\phi):=\phi(0)$ defines a distribution on $I$,
since
\[
 |\delta_0(\phi)|\leq\|\phi\|_{C^0(I)}.
\]
Since $I$ is compact, $\mathcal D(I)=C^\infty(I)$, with smoothness
up to the endpoints. Define
\[
 A\colon\mathcal D(I)\longrightarrow\mathcal D'(I),
 \qquad A\psi:=\psi(0)\delta_0.
\]
For each fixed $\phi$, $(A\psi)(\phi)=\psi(0)\phi(0)$ is continuous
in $\psi$; by the pointwise convergence topology, $A$ is continuous.
Its associated bilinear form has the joint bound
\[
 |B_A(\phi,\psi)|\leq\|\phi\|_{C^0(I)}\|\psi\|_{C^0(I)}.
\]
The product $I^2$ has a depth-two corner at $(0,0)$, and evaluation
\[
 K(\Phi):=\Phi(0,0),\qquad \Phi\in C^\infty(I^2),
\]
is a continuous distribution, since
$|K(\Phi)|\leq\|\Phi\|_{C^0(I^2)}$. For every pair of tests,
\[
 K(\phi\boxtimes\psi)=\phi(0)\psi(0)=(A\psi)(\phi).
\]
Uniqueness in Theorem~\ref{teo: nucleo Schwartz haces} shows that
$K=K_A=\delta_{(0,0)}$. Its support is exactly $\{(0,0)\}$: it vanishes
on tests supported away from that point, and every relative
neighborhood of the point contains a smooth cutoff equal to one
there, on which it has value one.

If tests were replaced by functions compactly supported in $(0,1)$
or $(0,1)^2$, both evaluations would vanish on every test.
That other choice of test space would lose the preceding operator.
The convention of tests smooth up to faces is therefore part of the
kernel formulation. It is compatible with defining Sobolev spaces
through weak identities in the interior: those identities do not,
by themselves, describe distributions supported on the boundary.
\end{example}

\begin{remark}[Dependencies of the kernel theorem with corners]
\label{obs:nucleo-esquinas-dependencias}
Extension of the kernel theorem to corners uses the product structure
of Chapter~1, partitions of unity, smooth orthant extension in
Appendix~B, nuclearity of the stages, and tensorization of test
spaces. Its proof uses neither Sobolev embeddings nor a trace theorem
on the total boundary. The results of
Section~\ref{sec:sobolev-haces-esquinas} permit subsequent study of
regularity and compactness for functions or sections on these
products, with the stated hypotheses. Existence of a distributional
kernel does not imply that it is a regular function or that the
operator satisfies elliptic boundary estimates.
\end{remark}

\begin{remark}[Variants of the kernel theorem] \label{rem:fuentes-teorema-nucleo-schwartz} The scalar Euclidean formulation can be found in \cite[Theorem 51.7]{Treves1967}. A proof by regularization and a formulation on manifolds appear in \cite[Theorem 3.1]{Ganglberger2010}. Taylor proves the compact-manifold version using Sobolev estimates and Hilbert--Schmidt operators in \cite[Chapter 4, §6]{TaylorPDE1}. For the geometric formulation with vector bundles over closed manifolds, see also \cite[Chapter 13, §3.5]{KrupkaSaunders2008}. The result here holds for noncompact manifolds and different bundles on the two factors. Control through the stages $K_j\times L_j$ replaces a global compactness hypothesis. \end{remark}
\chapter{Principal symbol and elliptic operators} \label{cap:simbolo-principal-operadores-elipticos}

Chapters~\ref{cap:operadores-diferenciales-parciales-haces} and~\ref{cap:distribuciones-en-haces-vectoriales} construct differential operators between sections, their formal adjoints, and their distributional extensions. Building on this foundation, this chapter isolates the highest-order part, formulates ellipticity intrinsically, and derives the regularity estimates that will be used to study the heat semigroup and fractional Sobolev spaces.

The principal part of an operator contains less information than the full operator, but it governs several of its fundamental analytic properties. On a compact manifold, invertibility of the principal symbol allows us to recover all derivatives of a section from the operator applied to that section and a zeroth-order term. Elliptic regularity and the equivalence of the Sobolev norm and the graph norm follow from this observation. The same symbol will later be the geometric object that enters index theory.

Throughout this chapter, $M$ will be a smooth manifold of dimension $n$, without boundary unless otherwise specified, and $\mathbf{E},\mathbf{F},\mathbf{G}\to M$ will be smooth real or complex vector bundles of finite rank. Whenever $L^2(M,\mathbf E)$ or Sobolev norms are involved, we fix a Riemannian metric $\mathbf{g}$ on $M$, bundle metrics $\mathbf{h}_{\mathbf{E}}$, $\mathbf{h}_{\mathbf{F}}$, and $\mathbf{h}_{\mathbf{G}}$ on the corresponding bundles, Hermitian in the complex case, and compatible connections $\nabla^{\mathbf E}$, $\nabla^{\mathbf F}$, and $\nabla^{\mathbf G}$. On tensor factors constructed from $TM$ and $T^*M$, we use the Levi--Civita connection of $\mathbf g$.

\begin{semblanzaHistorica}{Ellipticity: a local condition with global consequences} Elliptic theory took its modern form with the recognition that the highest-order terms control regularity, whereas lower-order terms can be treated as perturbations. The principal symbol encodes this dominant part on the cotangent bundle. Its invertibility away from the zero section is a pointwise condition, yet it leads to global estimates, finite-dimensional kernels, and stability of the index. This contrast between a microlocal hypothesis and global consequences will recur throughout the remainder of the book. \end{semblanzaHistorica}

\section{The principal symbol}

Recall that $\mathbf{PDO}^{(m)}(\mathbf{E},\mathbf{F})$ denotes the space of partial differential operators of order at most $m$ introduced in Definition~\ref{def:operadores-diferenciales-en-haces-pdo}. For $f\in C^\infty(M)$, we write $\operatorname{ad}(f)(P)=[P,f]$, where $[P,f](u)=P(fu)-fP(u)$. Proposition~\ref{prop: localidad respecto al diferencial del adjunto de Lie} shows that the pointwise evaluation of an iterated commutator of length $m$ depends only on the differentials of the functions involved.

\begin{definition}[Principal symbol] \label{def:simbolo-principal-pdo-haces} \index{principal symbol@principal symbol} Let $M$ be a smooth manifold without boundary, let $\mathbf{E},\mathbf{F}\to M$ be smooth vector bundles of finite rank, and let $P\in\mathbf{PDO}^{(m)}(\mathbf{E},\mathbf{F})$, with $m\geq1$. For $x\in M$, $\xi_1,\dots,\xi_m\in T_x^*M$, and $v\in \mathbf{E}_x$, choose functions $f_1,\dots,f_m\in C^\infty(M)$ such that $df_j(x)=\xi_j$ and define \[
\widetilde{\boldsymbol{\sigma}}_m(P)_x(\xi_1,\dots,\xi_m)v
:=
\frac{1}{m!}
\operatorname{ad}(f_1)\cdots\operatorname{ad}(f_m)(P)|_x(v).
\] The \textbf{polarized form of the principal symbol} is the map $\widetilde{\boldsymbol{\sigma}}_m(P)_x\colon (T_x^*M)^m\longrightarrow\operatorname{Hom}(\mathbf{E}_x,\mathbf{F}_x)$, and the \textbf{principal symbol} of $P$ is the map \[
\boldsymbol{\sigma}_m(P)\colon T^*M\longrightarrow\operatorname{Hom}(\mathbf{E},\mathbf{F}),
\qquad
\boldsymbol{\sigma}_m(P)(x,\xi)
:=
\widetilde{\boldsymbol{\sigma}}_m(P)_x(\xi,\dots,\xi).
\] If $m=0$, the operator $P$ is $C^\infty(M)$-linear and corresponds to a bundle homomorphism; in this case we set $\boldsymbol{\sigma}_0(P)(x):=P_x$. \end{definition}

The normalization by $m!$ ensures that, in coordinates, the symbol is exactly the homogeneous polynomial formed from the coefficients of order $m$.

\begin{proposition}[Well-definedness and elementary properties] \label{prop:simbolo-principal-bien-definido} Let $M$ be a smooth manifold without boundary, let $\mathbf{E},\mathbf{F}\to M$ be smooth vector bundles of finite rank, and let $P\in\mathbf{PDO}^{(m)}(\mathbf{E},\mathbf{F})$. The form $\widetilde{\boldsymbol{\sigma}}_m(P)_x$ is well defined, symmetric, and multilinear in the covariables. Moreover, $\boldsymbol{\sigma}_m(P)$ is smooth and satisfies \[
\boldsymbol{\sigma}_m(P)(x,t\xi)=t^m\boldsymbol{\sigma}_m(P)(x,\xi)
\] for all $x\in M$, $\xi\in T_x^*M$, and $t\in\mathbb R$. \end{proposition}

\begin{proof} Independence of the chosen functions is precisely the content of Proposition~\ref{prop: localidad respecto al diferencial del adjunto de Lie}. We prove linearity in the first covariable; the other arguments are treated in the same way. Let $\xi,\eta,\xi_2,\dots,\xi_m\in T_x^*M$ and $a,b\in\mathbb R$. Choose smooth functions $f,g,f_2,\dots,f_m$ such that $df(x)=\xi$, $dg(x)=\eta$, and $df_j(x)=\xi_j$. Since $d(af+bg)(x)=a\xi+b\eta$, the cited proposition allows us to use $af+bg$ to represent this covariable. Linearity of $f\mapsto\operatorname{ad}(f)$ then gives \[
\widetilde{\boldsymbol{\sigma}}_m(P)_x(a\xi+b\eta,\xi_2,\dots,\xi_m)=
a\widetilde{\boldsymbol{\sigma}}_m(P)_x(\xi,\xi_2,\dots,\xi_m)
+b\widetilde{\boldsymbol{\sigma}}_m(P)_x(\eta,\xi_2,\dots,\xi_m).
\] Lemma~\ref{lem:operadores-diferenciales-en-haces-pdo} allows us to interchange two consecutive commutators; hence the form is symmetric, and the linearity just proved in the first entry implies linearity in every entry.

It remains to check smoothness. In a chart $(U,\varphi)$, with $\varphi=(x^1,\dots,x^n)$, and local frames $(\mathbf{e}_1,\dots,\mathbf{e}_r)$ of $\mathbf{E}$ and $(\mathbf{f}_1,\dots,\mathbf{f}_s)$ of $\mathbf{F}$, Corollary~\ref{cor:pdo-expresion-con-marcos} allows us to write \[
P\left(\displaystyle\sum_{a=1}^r u^a\mathbf{e}_a\right)
=
\displaystyle\sum_{b=1}^s
\left(
\displaystyle\sum_{|\alpha|\leq m}\displaystyle\sum_{a=1}^r
(A_\alpha)_{ba}D^\alpha u^a
\right)\mathbf{f}_b,
\] where $A_\alpha\colon U\longrightarrow M_{s\times r}(\mathbb K)$ is smooth. For $x\in U$ and $\xi=\displaystyle\sum_{i=1}^n\xi_i\,\mathbf{d}x^i|_x$, consider near $x$ the linear coordinate function $\ell_\xi:=\displaystyle\sum_{i=1}^n\xi_i x^i$. Multiplying it by a smooth function equal to one on a neighborhood of $x$ produces a function in $C^\infty(M)$ with the same differential at $x$. Proposition~\ref{prop: localidad respecto al diferencial del adjunto de Lie} therefore allows us to use $\ell_\xi$ in the $m$ entries defining $\boldsymbol{\sigma}_m(P)(x,\xi)$.

Since $[\partial_i,\ell_\xi]=\xi_i\operatorname{Id}$ and $[AB,\ell_\xi]=A[B,\ell_\xi]+[A,\ell_\xi]B$, each commutator replaces one of the $|\alpha|$ differential factors of $D^\alpha$ with the corresponding component of $\xi$. If $|\alpha|<m$, the result vanishes after $m$ commutators. If $|\alpha|=m$, a term survives only if each of the $m$ factors is replaced exactly once, and the commutators can select these factors in $m!$ orders. Consequently, \[
\operatorname{ad}(\ell_\xi)^m(D^\alpha)
=
\begin{cases}
m!\,\xi^\alpha\operatorname{Id},& |\alpha|=m,\\
0,&|\alpha|<m.
\end{cases}
\] The matrices $A_\alpha$ act by multiplication and commute with scalar functions. Dividing by the normalization $m!$ in the definition yields \[
\boldsymbol{\sigma}_m(P)(x,\xi)
=
\displaystyle\sum_{|\alpha|=m}A_\alpha(x)\xi^\alpha.
\] The coefficients $A_\alpha$ are smooth, and the right-hand side is a homogeneous polynomial of degree $m$ in $\xi$; therefore $\boldsymbol{\sigma}_m(P)$ is smooth and $\boldsymbol{\sigma}_m(P)(x,t\xi)=t^m\boldsymbol{\sigma}_m(P)(x,\xi)$. \end{proof}

\begin{proposition}[Local expression] \label{prop:formula-local-simbolo-principal} Let $M$ be a smooth manifold without boundary, let $\mathbf{E},\mathbf{F}\to M$ be smooth vector bundles of finite rank, and let $P\in\mathbf{PDO}^{(m)}(\mathbf{E},\mathbf{F})$. In a chart and local frames in which $\displaystyle P=\displaystyle\sum_{|\alpha|\leq m}A_\alpha D^\alpha$, if $\xi=\displaystyle\sum_{j=1}^n\xi_jdx^j|_x$, then \[
\boldsymbol{\sigma}_m(P)(x,\xi)
=
\displaystyle\sum_{|\alpha|=m}A_\alpha(x)\xi^\alpha,
\qquad
\xi^\alpha:=\prod_{j=1}^n\xi_j^{\alpha_j}.
\] Consequently, $P$ has order at most $m-1$ if and only if $\boldsymbol{\sigma}_m(P)=0$. \end{proposition}

\begin{proof} This formula is the diagonal computation carried out in the proof of Proposition~\ref{prop:simbolo-principal-bien-definido}. If $P$ has order at most $m-1$, no coefficients $A_\alpha$ with $|\alpha|=m$ occur, and the symbol vanishes. Conversely, if $\boldsymbol{\sigma}_m(P)=0$, every entry of the polynomial matrix $\displaystyle\sum_{|\alpha|=m}A_\alpha(x)\xi^\alpha$ vanishes for all $\xi\in\mathbb R^n$. Uniqueness of polynomial coefficients implies $A_\alpha(x)=0$ for every $|\alpha|=m$ and every $x$ in the chart. Corollary~\ref{cor:pdo-expresion-con-marcos} then shows that $P\in\mathbf{PDO}^{(m-1)}(\mathbf{E},\mathbf{F})$. \end{proof}

\begin{remark} The preceding definition does not include powers of $i$. With this convention, the symbol of $-\partial_j^2$ is $-\xi_j^2$. Applying the Fourier transform instead introduces the multiplier $(i\xi)^\alpha$. We will explicitly account for this distinction in the proof of the elliptic estimates. \end{remark}

\subsection{Elementary symbol calculus}

The principal symbol records the highest-order part of an operator. Before using it to formulate ellipticity, we need to check how it behaves under composition and passage to adjoints.

\begin{proposition}[Symbol of a composition] \label{prop:simbolo-composicion-pdo} Let $M$ be a smooth manifold without boundary, and let $\mathbf{E},\mathbf{F},\mathbf{G}\to M$ be smooth vector bundles of finite rank. Let $P\in\mathbf{PDO}^{(m)}(\mathbf{E},\mathbf{F})$ and $Q\in\mathbf{PDO}^{(\ell)}(\mathbf{F},\mathbf{G})$. Then $Q\circ P\in\mathbf{PDO}^{(m+\ell)}(\mathbf{E},\mathbf{G})$ and \[
\boldsymbol{\sigma}_{m+\ell}(Q\circ P)(x,\xi)
=
\boldsymbol{\sigma}_\ell(Q)(x,\xi)\circ\boldsymbol{\sigma}_m(P)(x,\xi).
\] \end{proposition}

\begin{proof} Membership in $\mathbf{PDO}^{(m+\ell)}(\mathbf{E},\mathbf{G})$ was proved in Proposition~\ref{prop: composicion de pdos es pdo}. In a chart and local frames, write \[
P=\displaystyle\sum_{|\alpha|\leq m}A_\alpha D^\alpha,
\qquad
Q=\displaystyle\sum_{|\beta|\leq\ell}B_\beta D^\beta.
\] Upon expanding $Q(Pu)$, terms of order $m+\ell$ arise only when all the derivatives $D^\beta$ act on $D^\alpha u$ and none acts on $A_\alpha$. Thus the highest-order part is \[
\displaystyle\sum_{|\beta|=\ell}\displaystyle\sum_{|\alpha|=m}
B_\beta A_\alpha D^{\alpha+\beta}u.
\] Evaluating the symbol at $\xi$ gives \[
\displaystyle\sum_{|\beta|=\ell}\displaystyle\sum_{|\alpha|=m}
B_\beta(x)A_\alpha(x)\xi^{\alpha+\beta}
=
\left(\displaystyle\sum_{|\beta|=\ell}B_\beta(x)\xi^\beta\right)
\left(\displaystyle\sum_{|\alpha|=m}A_\alpha(x)\xi^\alpha\right),
\] which is the asserted identity. \end{proof}

The formal Hermitian adjoint was already constructed in Theorem~\ref{cor: adjunto global en coordenadas locales (existencia y expresion)}. We need only record the effect of that construction on the principal part.

\begin{proposition}[Symbol of the formal Hermitian adjoint] \label{prop:simbolo-adjunto-formal} Let $(M,\mathbf{g})$ be a Riemannian manifold without boundary, and let $\mathbf{E},\mathbf{F}\to M$ be smooth vector bundles of finite rank equipped with real or Hermitian bundle metrics $\mathbf{h}_{\mathbf{E}}$ and $\mathbf{h}_{\mathbf{F}}$. If $P\in\mathbf{PDO}^{(m)}(\mathbf{E},\mathbf{F})$ and $P_h^*$ is its formal Hermitian adjoint, then \[
\boldsymbol{\sigma}_m(P_h^*)(x,\xi)
=
(-1)^m\boldsymbol{\sigma}_m(P)(x,\xi)^*,
\] where the adjoint on the right is taken between the fibers $(\mathbf{E}_x,\mathbf{h}_{\mathbf{E}})$ and $(\mathbf{F}_x,\mathbf{h}_{\mathbf{F}})$. \end{proposition}

\begin{proof} In orthonormal or unitary local frames, write \[
P\left(\sum_{i=1}^r u^i\mathbf e_i\right)
=
\sum_{j=1}^s
\left(
\sum_{i=1}^r\sum_{|\alpha|\leq m}
a_{\alpha i}^{\,j}D_\phi^\alpha u^i
\right)\mathbf f_j.
\] Theorem~\ref{cor: adjunto global en coordenadas locales (existencia y expresion)} shows that the component along $\mathbf e_i$ of $\displaystyle P_h^*(\displaystyle\sum_{j=1}^{s}v^j\mathbf f_j)$ is \[
\frac{1}{\sqrt{\det(\mathbf{g})}}
\sum_{j=1}^s\sum_{|\alpha|\leq m}(-1)^{|\alpha|}
D_\phi^\alpha\!\left(
\sqrt{\det(\mathbf{g})}\,
\overline{a_{\alpha i}^{\,j}}\,v^j
\right).
\] Terms of order $m$ arise only when all $|\alpha|$ derivatives act on $v^j$. If \[
A_\alpha=(a_{\alpha i}^{\,j})\colon
U\longrightarrow M_{s\times r}(\mathbb K),
\] their coefficients form, pointwise, the matrix $(-1)^mA_\alpha^*$ for $|\alpha|=m$. Proposition~\ref{prop:formula-local-simbolo-principal} gives \[
\boldsymbol{\sigma}_m(P_h^*)(x,\xi)
=
(-1)^m\displaystyle\sum_{|\alpha|=m}A_\alpha(x)^*\xi^\alpha
=
(-1)^m\boldsymbol{\sigma}_m(P)(x,\xi)^*.
\] \end{proof}

\begin{corollary}[Symbol under conjugation] \label{cor:simbolo-principal-conjugacion-operadores} Let $M$ be a smooth manifold without boundary, let $\mathbf{E},\mathbf{F}\to M$ be smooth complex vector bundles of finite rank, and let $P\in\mathbf{PDO}^{(m)}(\mathbf{E},\mathbf{F})$. Then \[
 \boldsymbol{\sigma}_m(\overline P)(x,\xi)
 =\overline{\boldsymbol{\sigma}_m(P)(x,\xi)}
 \colon\overline{\mathbf{E}_x}\longrightarrow\overline{\mathbf{F}_x}.
\] In particular, $P$ is elliptic if and only if $\overline P$ is elliptic. \end{corollary}

\begin{proof} In conjugate frames, the principal coefficients of $\overline P$ are $\overline{A_\alpha}$ by Proposition~\ref{prop:conjugacion-operadores-diferenciales-haces}. Proposition~\ref{prop:formula-local-simbolo-principal} gives the formula; conjugation preserves invertibility. \end{proof}

\begin{corollary}[Symbol of the formal transpose] \label{cor:simbolo-transpuesto-formal} Let $(M,\mathbf{g})$ be a Riemannian manifold without boundary, let $\mathbf{E},\mathbf{F}\to M$ be smooth vector bundles of finite rank, and let $P\in\mathbf{PDO}^{(m)}(\mathbf{E},\mathbf{F})$. If $P'$ is its bilinear formal transpose, then \begin{equation}
\label{eq:simbolo-transpuesto-formal}
 \boldsymbol{\sigma}_m(P')(x,\xi)
 =(-1)^m\boldsymbol{\sigma}_m(P)(x,\xi)^T
 \colon \mathbf{F}_x^*\longrightarrow \mathbf{E}_x^*,
\end{equation} where $T$ is the algebraic transpose, without conjugation. Thus $P$ is elliptic if and only if $P'$ is elliptic. \end{corollary}

\begin{proof} Corollary~\ref{cor:formula-local-transpuesto-formal} shows that the coefficient of $D^\alpha$ of order $m$ is $(-1)^mA_\alpha^T$; every derivative falling on the coordinate weight of the measure or on $A_\alpha$ has lower order. The local formula for the symbol yields \eqref{eq:simbolo-transpuesto-formal}. The transpose of a finite-dimensional map is invertible exactly when the map itself is invertible. \end{proof}

\begin{example}[The connection] \label{ej:simbolo-conexion} Let $M$ be a smooth manifold without boundary, and let $\mathbf{E}\to M$ be a smooth vector bundle with connection $\nabla^{\mathbf{E}}$. By Proposition~\ref{prop: conexion es pdo de orden 1}, $\nabla^{\mathbf{E}}$ is a first-order differential operator from $\mathbf{E}$ to $T^*M\otimes \mathbf{E}$. The Leibniz rule gives \[
[\nabla^{\mathbf{E}},f]\mathbf{u}=df\otimes \mathbf{u},
\] and consequently \[
\boldsymbol{\sigma}_1(\nabla^{\mathbf{E}})(x,\xi)v=\xi\otimes v.
\] This symbol is injective for $\xi\neq0$, although it is not surjective when $n>1$. \end{example}

\begin{example}[The Bochner Laplacian] \label{ej:simbolo-laplaciano-bochner} Let $(M,\mathbf{g})$ be a Riemannian manifold without boundary, and let $\mathbf{E}\to M$ be a smooth vector bundle with a real or Hermitian bundle metric $\mathbf{h}_{\mathbf{E}}$ and a compatible connection $\nabla^{\mathbf{E}}$. The Bochner Laplacian $\Delta_B=(\nabla^{\mathbf{E}})_h^*\nabla^{\mathbf{E}}$, defined in Definition~\ref{def:operadores-diferenciales-en-haces-laplaciano-de-bochner}, satisfies \[
\boldsymbol{\sigma}_2(\Delta_B)(x,\xi)
=
-\|\xi\|_{\mathbf{g}}^2\operatorname{Id}_{\mathbf{E}_x}.
\] This identity follows from Propositions~\ref{prop:simbolo-composicion-pdo} and~\ref{prop:simbolo-adjunto-formal}, since \[
\boldsymbol{\sigma}_1((\nabla^{\mathbf{E}})_h^*)(x,\xi)(\eta\otimes v)
=-\langle\xi,\eta\rangle_{\mathbf{g}}v.
\] \end{example}

\section{Elliptic operators}

Ellipticity expresses the invertibility of the principal symbol in every nonzero cotangent direction. This condition allows derivatives of the solution to be recovered from the equation and leads to regularity estimates.

\begin{definition}[Ellipticity] \label{def:operador-eliptico-haces} \index{elliptic operator@elliptic operator} Let $M$ be a smooth manifold without boundary, let $\mathbf{E},\mathbf{F}\to M$ be smooth vector bundles of finite rank, and let $P\in\mathbf{PDO}^{(m)}(\mathbf{E},\mathbf{F})$, with $m\geq1$. We say that $P$ is \textbf{elliptic} if, for every $x\in M$ and every $\xi\in T_x^*M\setminus\{0\}$, the map \[
\boldsymbol{\sigma}_m(P)(x,\xi)\colon \mathbf{E}_x\longrightarrow \mathbf{F}_x
\] is a linear isomorphism. \end{definition}

The definition requires $\mathbf{E}$ and $\mathbf{F}$ to have the same rank. When the symbol is only injective or only surjective, one speaks of an overdetermined or underdetermined elliptic operator, respectively. These variants will be useful in index theory, but are not needed for the estimates in this chapter.

\begin{proposition} \label{prop:elipticidad-adjunto-composicion} Let $(M,\mathbf{g})$ be a Riemannian manifold without boundary, let $\mathbf{E},\mathbf{F},\mathbf{G}\to M$ be smooth vector bundles of finite rank equipped with bundle metrics, and let $P\in\mathbf{PDO}^{(m)}(\mathbf{E},\mathbf{F})$. \begin{enumerate}[label=(\alph*)] \item $P$ is elliptic if and only if $P_h^*$ is elliptic; \item if $P$ is elliptic, then $P_h^*P$ and $PP_h^*$ are elliptic of order $2m$; \item if $Q\in\mathbf{PDO}^{(\ell)}(\mathbf{F},\mathbf{G})$ and both $P$ and $Q$ are elliptic, then $Q\circ P$ is elliptic. \end{enumerate} \end{proposition}

\begin{proof} Part (a) follows from Proposition~\ref{prop:simbolo-adjunto-formal}, because a linear map between finite-dimensional spaces is invertible if and only if its adjoint is invertible. For (b), Propositions~\ref{prop:simbolo-composicion-pdo} and~\ref{prop:simbolo-adjunto-formal} give \[
\boldsymbol{\sigma}_{2m}(P_h^*P)(x,\xi)
=
(-1)^m\boldsymbol{\sigma}_m(P)(x,\xi)^*\boldsymbol{\sigma}_m(P)(x,\xi),
\] which is invertible for $\xi\neq0$. The argument for $PP_h^*$ is identical. Part (c) follows because a composition of isomorphisms is an isomorphism. \end{proof}

\begin{example} The Bochner Laplacian is elliptic by Example~\ref{ej:simbolo-laplaciano-bochner}. The Weitzenböck identity, Theorem~\ref{teo:operadores-diferenciales-en-haces-identidad-de-weitzenbock}, shows that the bundle-valued Hodge Laplacian differs from the Bochner Laplacian by a zeroth-order endomorphism. They therefore have the same principal symbol and are both elliptic. \end{example}

\begin{proposition}[Differential symbol of a Dirac-type operator] \label{prop:simbolo-elipticidad-dirac} \index{Dirac operator!principal symbol} Let \((\mathbf S,c,\nabla^{\mathbf S})\) be a Clifford module with connection, and let \(D=c\circ\nabla^{\mathbf S}\) be the operator in Definition~\ref{def:operador-tipo-dirac}. Its differential principal symbol is \begin{equation}
 \boldsymbol\sigma_1(D)(x,\xi)=c(\xi)\colon
 \mathbf S_x\longrightarrow\mathbf S_x.
\label{eq:simbolo-diferencial-dirac}
\end{equation} In particular, \(D\) is elliptic and, for \(\xi\neq0\), \[
 \boldsymbol\sigma_1(D)(x,\xi)^{-1}
 =-\frac{c(\xi)}{|\xi|_{\mathbf g}^2}.
\] If the module is graded, the block \(D^+\colon\Gamma(\mathbf S^+)\to\Gamma(\mathbf S^-)\) is elliptic and its symbol is the restriction \(c(\xi)\colon\mathbf S_x^+\to\mathbf S_x^-\). \end{proposition}

\begin{proof} By \eqref{eq:conmutador-dirac-funcion}, for any function \(f\) with \(df_x=\xi\), \[
 \boldsymbol\sigma_1(D)(x,\xi)
 =[D,f]_x=c(df_x)=c(\xi).
\] The Clifford relation \(c(\xi)^2=-|\xi|_{\mathbf g}^2I\) gives the stated inverse. In the graded case, \(c(\xi)\) is odd and its square satisfies the same identity on each summand, so both restrictions are isomorphisms for \(\xi\neq0\). \end{proof}

\section{Elliptic estimates on vector bundles}

We first prove the estimate for systems with constant coefficients, then freeze the coefficients at a point, and finally globalize the result by means of a partition of unity. The only ingredient from Fourier analysis is Plancherel's identity.

\subsection{An interpolation inequality}

Local estimates involve lower-order terms that must be absorbed. The following inequality quantifies this absorption in terms of a small portion of the higher-order norm and a lower-order norm.

\begin{lemma}[Interpolation between integer orders] \label{lem:interpolacion-sobolev-orden-inferior} Let $m\in\mathbb N$ and $0\leq j<m$. For every $\varepsilon>0$, there exists $C_\varepsilon>0$ such that \[
\|u\|_{H^j(\mathbb R^n,\mathbb K^r)}
\leq
\varepsilon\|u\|_{H^m(\mathbb R^n,\mathbb K^r)}
+C_\varepsilon\|u\|_{L^2(\mathbb R^n,\mathbb K^r)}
\] for every $u\in H^m(\mathbb R^n,\mathbb K^r)$. \end{lemma}

\begin{proof} We use the Fourier transform with unitary normalization, write $\langle\xi\rangle:=(1+\|\xi\|^2)^{\displaystyle\frac{1}{2}}$, and set \[
\|u\|_{s,\mathbf{F}}:=
\|\langle\xi\rangle^s\widehat u\|_{L^2(\mathbb R^n,\mathbb K^r)}.
\] For each $s\in\mathbb N_0$, the Leibniz formula for $(1+\|\xi\|^2)^s$ and Plancherel's identity give constants $a_s,b_s>0$, depending only on $n$, $r$, $s$, and the chosen normalization of the Sobolev norm, such that \[
a_s\|u\|_{H^s(\mathbb R^n,\mathbb K^r)}
\leq\|u\|_{s,\mathbf{F}}
\leq b_s\|u\|_{H^s(\mathbb R^n,\mathbb K^r)}.
\]

Let $R\geq1$ and set $\mathbf{E}_R:=\{\xi\in\mathbb R^n\mid\langle\xi\rangle\leq R\}$. On $\mathbf{E}_R$ we have $\langle\xi\rangle^j\leq R^j$, whereas on $\mathbf{E}_R^c$ we have $\langle\xi\rangle^{j-m}\leq R^{j-m}$. By the triangle inequality in $L^2(\mathbb R^n,\mathbb K^r)$, applied componentwise, and Plancherel's identity, \[
\begin{aligned}
\|u\|_{j,\mathbf{F}}
&=\|\langle\xi\rangle^j\widehat u\|_{L^2(\mathbb R^n,\mathbb K^r)}\\
&\leq
\|\mathbf 1_{\mathbf{E}_R}\langle\xi\rangle^j\widehat u\|_{L^2(\mathbb R^n,\mathbb K^r)}
+\|\mathbf 1_{\mathbf{E}_R^c}\langle\xi\rangle^j\widehat u\|_{L^2(\mathbb R^n,\mathbb K^r)}\\
&\leq
R^j\|\widehat u\|_{L^2(\mathbb R^n,\mathbb K^r)}
+R^{j-m}\|\langle\xi\rangle^m\widehat u\|_{L^2(\mathbb R^n,\mathbb K^r)}\\
&=
R^j\|u\|_{L^2(\mathbb R^n,\mathbb K^r)}
+R^{j-m}\|u\|_{m,\mathbf{F}}.
\end{aligned}
\] Let $\delta:=\displaystyle\frac{a_j\varepsilon}{b_m}$ and choose $R_\varepsilon:=\displaystyle\max\left\{1,\delta^{-\displaystyle\frac{1}{m-j}}\right\}$. Then $R_\varepsilon^{j-m}\leq\delta$ and, using the two preceding inequalities, \[
\|u\|_{H^j(\mathbb R^n,\mathbb K^r)}
\leq a_j^{-1}\|u\|_{j,\mathbf{F}}
\leq\varepsilon\|u\|_{H^m(\mathbb R^n,\mathbb K^r)}
+a_j^{-1}R_\varepsilon^j\|u\|_{L^2(\mathbb R^n,\mathbb K^r)}.
\] Thus it suffices to take $C_\varepsilon:=a_j^{-1}R_\varepsilon^j$. \end{proof}

\subsection{Systems with constant coefficients}

Freezing the coefficients first reduces the problem to a translation-invariant system. In this model, the Fourier transform converts ellipticity into an algebraic estimate for the symbol.

\begin{lemma}[Estimate for a constant-coefficient elliptic system] \label{lem:estimacion-eliptica-coeficientes-constantes} Let $n,r,m\in\mathbb N$, let $\mathbb K\in\{\mathbb R,\mathbb C\}$, and let $A_\alpha\in\operatorname{End}(\mathbb K^r)$ for $|\alpha|\leq m$. Consider the operator \[
P_0(D)\colon C_c^\infty(\mathbb R^n,\mathbb K^r)
\longrightarrow C_c^\infty(\mathbb R^n,\mathbb K^r),
\qquad
P_0(D)u=\displaystyle\sum_{|\alpha|\leq m}A_\alpha D^\alpha u,
\] with constant coefficients, and suppose that \[
p_m(\xi):=\displaystyle\sum_{|\alpha|=m}A_\alpha\xi^\alpha
\] is invertible for $\xi\neq0$. There exists $C>0$ such that \[
\|u\|_{H^m(\mathbb R^n,\mathbb K^r)}
\leq
C\bigl(\|P_0(D)u\|_{L^2(\mathbb R^n,\mathbb K^r)}
+\|u\|_{L^2(\mathbb R^n,\mathbb K^r)}\bigr)
\] for every $u\in C_c^\infty(\mathbb R^n,\mathbb K^r)$. \end{lemma}

\begin{proof} If $\mathbb K=\mathbb R$, we complexify the system to use the Fourier transform; the resulting estimate is then restricted to real-valued functions. The function $(\xi,v)\mapsto\|p_m(\xi)v\|_{\mathbb C^r}$ is continuous on the compact set $S^{n-1}\times S^{2r-1}$ and does not vanish. Thus there exists $c>0$ such that \[
\|p_m(\xi)v\|_{\mathbb C^r}\geq c\|\xi\|^m\|v\|_{\mathbb C^r}
\] for all $\xi\in\mathbb R^n$ and $v\in\mathbb C^r$.

We must distinguish the preceding principal symbol from the Fourier multiplier. Set \[
p^\#(\xi):=\displaystyle\sum_{|\alpha|\leq m}A_\alpha(i\xi)^\alpha,
\qquad
p_m^\#(\xi):=\displaystyle\sum_{|\alpha|=m}A_\alpha(i\xi)^\alpha=i^mp_m(\xi).
\] Since $|i^m|=1$, the preceding estimate gives $\|p_m^\#(\xi)v\|_{\mathbb C^r}\geq c\|\xi\|^m\|v\|_{\mathbb C^r}$. If \[
r^\#(\xi):=p^\#(\xi)-p_m^\#(\xi),
\] then there exists $C_0>0$ such that $\|r^\#(\xi)v\|_{\mathbb C^r}\leq C_0(1+\|\xi\|)^{m-1}\|v\|_{\mathbb C^r}$. Consequently, \[
\|p^\#(\xi)v\|_{\mathbb C^r}
\geq
c\|\xi\|^m\|v\|_{\mathbb C^r}-C_0(1+\|\xi\|)^{m-1}\|v\|_{\mathbb C^r}.
\] Choose $R>1$ so that $C_0(1+\|\xi\|)^{m-1}\leq\displaystyle\frac{c}{2}\|\xi\|^m$ for $\|\xi\|\geq R$. In this region, \[
(1+\|\xi\|^2)^{\frac{m}{2}}\|v\|_{\mathbb C^r}
\leq
C_1\|p^\#(\xi)v\|_{\mathbb C^r}.
\] For $\|\xi\|<R$, we have $(1+\|\xi\|^2)^{\displaystyle\frac{m}{2}}\|v\|_{\mathbb C^r}\leq C_R\|v\|_{\mathbb C^r}$. Therefore \[
(1+\|\xi\|^2)^{\frac{m}{2}}\|\widehat u(\xi)\|_{\mathbb C^r}
\leq
C\bigl(\|p^\#(\xi)\widehat u(\xi)\|_{\mathbb C^r}+\|\widehat u(\xi)\|_{\mathbb C^r}\bigr).
\] We square, integrate, and apply Plancherel's identity. Since \[
\widehat{P_0(D)u}(\xi)=p^\#(\xi)\widehat u(\xi),
\] the asserted estimate follows. \end{proof}

\begin{lemma}[Constant-coefficient estimate at an arbitrary Sobolev order] \label{lem:estimacion-eliptica-coeficientes-constantes-orden-arbitrario} Let $n,r,m\in\mathbb N$, let $\mathbb K\in\{\mathbb R,\mathbb C\}$, and let $A_\alpha\in\operatorname{End}(\mathbb K^r)$ for $|\alpha|=m$. Let the operator \[
P_m(D)\colon\mathcal S(\mathbb R^n,\mathbb K^r)
\longrightarrow\mathcal S(\mathbb R^n,\mathbb K^r),
\qquad
P_m(D)v=\displaystyle\sum_{|\alpha|=m}A_\alpha D^\alpha v,
\] be a homogeneous system with constant coefficients, and suppose that its symbol \[
p_m(\xi):=\displaystyle\sum_{|\alpha|=m}A_\alpha\xi^\alpha
\] is invertible for $\xi\neq0$. For every $s\in\mathbb R$, there exists $C_s>0$ such that \[
\|v\|_{H^{s+m}(\mathbb R^n,\mathbb K^r)}
\leq
C_s\bigl(
\|P_m(D)v\|_{H^s(\mathbb R^n,\mathbb K^r)}
+
\|v\|_{H^{s+m-1}(\mathbb R^n,\mathbb K^r)}
\bigr)
\] for every $v\in\mathcal S(\mathbb R^n,\mathbb K^r)$. \end{lemma}

\begin{proof} As in the proof of Lemma~\ref{lem:estimacion-eliptica-coeficientes-constantes}, ellipticity and compactness of the sphere give $c>0$ such that \[
\|p_m^\#(\xi)z\|_{\mathbb C^r}
\geq
c\|\xi\|^m\|z\|_{\mathbb C^r},
\qquad
p_m^\#(\xi):=\displaystyle\sum_{|\alpha|=m}A_\alpha(i\xi)^\alpha,
\] for all $\xi\in\mathbb R^n$ and $z\in\mathbb C^r$. Choose $R\geq1$. If $\|\xi\|\geq R$, the preceding inequality implies \[
\langle\xi\rangle^{s+m}\|z\|_{\mathbb C^r}
\leq
C_{s,R}\langle\xi\rangle^s\|p_m^\#(\xi)z\|_{\mathbb C^r}.
\] If $\|\xi\|<R$, then \[
\langle\xi\rangle^{s+m}\|z\|_{\mathbb C^r}
\leq
C_{s,R}\langle\xi\rangle^{s+m-1}\|z\|_{\mathbb C^r}.
\] Combining the two regions gives, for every $\xi$ and every $z$, \[
\langle\xi\rangle^{s+m}\|z\|_{\mathbb C^r}
\leq
C_s\left(
\langle\xi\rangle^s\|p_m^\#(\xi)z\|_{\mathbb C^r}
+
\langle\xi\rangle^{s+m-1}\|z\|_{\mathbb C^r}
\right).
\] We apply this inequality to $z=\widehat v(\xi)$, square, and integrate. Since $\widehat{P_m(D)v}=p_m^\#\widehat v$, Plancherel's identity completes the proof. \end{proof}

\subsection{Freezing coefficients}

On a small ball, a variable-coefficient operator is compared with the constant-coefficient operator obtained by fixing at the center only the coefficients of its principal part. This part determines ellipticity. The principal error retains order $m$, but its coefficients become small; the remaining terms need not be small because they have order at most $m-1$.

\begin{lemma}[Local estimate for elliptic systems] \label{lem:estimacion-eliptica-local-sistemas} Let $n,r,m\in\mathbb N$, let $\mathbb K\in\{\mathbb R,\mathbb C\}$, let $U\subseteq\mathbb R^n$ be open, and let $A_\alpha\in C^\infty(U;\operatorname{End}(\mathbb K^r))$. Let \[
P\colon C_c^\infty(U,\mathbb K^r)
\longrightarrow C_c^\infty(U,\mathbb K^r),
\qquad
Pu=\displaystyle\sum_{|\alpha|\leq m}A_\alpha(x)D^\alpha u,
\] be an elliptic system with smooth coefficients. For every $x_0\in U$, there exist an open subset $V$, with $x_0\in V\Subset U$, and a constant $C>0$ such that \[
\|u\|_{H^m(\mathbb R^n,\mathbb K^r)}
\leq
C\bigl(\|Pu\|_{L^2(\mathbb R^n,\mathbb K^r)}
+\|u\|_{L^2(\mathbb R^n,\mathbb K^r)}\bigr)
\] for every $u\in C_c^\infty(V,\mathbb K^r)$, extended by zero outside $V$. \end{lemma}

\begin{proof} We freeze only the principal coefficients at $x_0$ and write \[
P=P_{x_0}^{(m)}+R_1+R_0,
\] where \[
P_{x_0}^{(m)}:=\displaystyle\sum_{|\alpha|=m}A_\alpha(x_0)D^\alpha,
\qquad
R_1:=\displaystyle\sum_{|\alpha|=m}\bigl(A_\alpha(x)-A_\alpha(x_0)\bigr)D^\alpha,
\] and $R_0$ collects the terms of order at most $m-1$. For every $W\Subset U$, the operators $P_{x_0}^{(m)}$, $R_1$, and $R_0$ send $C_c^\infty(W,\mathbb K^r)$ into $C_c^\infty(U,\mathbb K^r)$. Lemma~\ref{lem:estimacion-eliptica-coeficientes-constantes} applied to $P_{x_0}^{(m)}$ gives \[
\|u\|_{H^m(\mathbb R^n,\mathbb K^r)}
\leq
C_0\bigl(\|P_{x_0}^{(m)}u\|_{L^2(\mathbb R^n,\mathbb K^r)}+\|u\|_{L^2(\mathbb R^n,\mathbb K^r)}\bigr).
\] Choose $\rho_0>0$ so that $\overline{B(x_0,\rho_0)}\subseteq U$. For $|\alpha|=m$, set \[
C_\alpha(\rho_0)
:=
\sup_{z\in\overline{B(x_0,\rho_0)}}
\sum_{j=1}^n
\left\|\frac{\partial A_\alpha}{\partial x^j}(z)\right\|,
\qquad
C_*(\rho_0):=\sum_{|\alpha|=m}C_\alpha(\rho_0).
\] Hadamard's lemma~\ref{lem:hadamard-finito-dimensional}, applied to the matrix entries of $A_\alpha$, gives, for $x\in B(x_0,\rho_0)$, \begin{equation}
A_\alpha(x)-A_\alpha(x_0)
=
\sum_{j=1}^n(x^j-x_0^j)
\int_0^1
\frac{\partial A_\alpha}{\partial x^j}
\bigl(x_0+t(x-x_0)\bigr)\,dt.
\label{eq:variacion-coeficientes-principales-hadamard}
\end{equation} Since the segment lies in the ball, for $0<\rho\leq \rho_0$ we obtain the quantitative bound \begin{equation}
\sup_{x\in B(x_0,\rho)}
\|A_\alpha(x)-A_\alpha(x_0)\|
\leq \rho C_\alpha(\rho_0).
\label{eq:cota-coeficiente-principal-radio}
\end{equation} Choose $\rho\in(0,\rho_0)$ small enough that \[
\rho C_*(\rho_0)\leq\frac{1}{4C_0C_1},
\] and take $V=B(x_0,\rho)$. If $C_*(\rho_0)=0$, any $\rho<\rho_0$ satisfies this condition. Then \[
\sup_{x\in V}\displaystyle\sum_{|\alpha|=m}\|A_\alpha(x)-A_\alpha(x_0)\|
\leq
\frac{1}{4C_0C_1},
\] where $C_1$ satisfies $\displaystyle\sum_{|\alpha|=m}\|D^\alpha u\|_{L^2(\mathbb R^n,\mathbb K^r)}\leq C_1\|u\|_{H^m(\mathbb R^n,\mathbb K^r)}$. Thus \[
\|R_1u\|_{L^2(\mathbb R^n,\mathbb K^r)}
\leq
\frac{1}{4C_0}\|u\|_{H^m(\mathbb R^n,\mathbb K^r)}.
\] Since $R_0$ has order at most $m-1$, there exists $C_2>0$ such that $\|R_0u\|_{L^2(\mathbb R^n,\mathbb K^r)}\leq C_2\|u\|_{H^{m-1}(\mathbb R^n,\mathbb K^r)}$. Applying Lemma~\ref{lem:interpolacion-sobolev-orden-inferior} with $\varepsilon=(4C_0C_2)^{-1}$ yields \[
\|R_0u\|_{L^2(\mathbb R^n,\mathbb K^r)}
\leq
\frac{1}{4C_0}\|u\|_{H^m(\mathbb R^n,\mathbb K^r)}+C_3\|u\|_{L^2(\mathbb R^n,\mathbb K^r)}.
\] Since $P_{x_0}^{(m)}u=Pu-R_1u-R_0u$, \[
\|u\|_{H^m(\mathbb R^n,\mathbb K^r)}
\leq
C_0\|Pu\|_{L^2(\mathbb R^n,\mathbb K^r)}
+\frac{1}{2}\|u\|_{H^m(\mathbb R^n,\mathbb K^r)}
+C_4\|u\|_{L^2(\mathbb R^n,\mathbb K^r)}.
\] Moving the middle term to the left-hand side completes the proof. \end{proof}

\begin{proposition}[Local estimate at an arbitrary Sobolev order] \label{prop:estimacion-eliptica-local-orden-arbitrario} Let $n,r,m\in\mathbb N$, let $\mathbb K\in\{\mathbb R,\mathbb C\}$, let $U\subseteq\mathbb R^n$ be open, and let $A_\alpha\in C^\infty(U;\operatorname{End}(\mathbb K^r))$. Let \[
P\colon C_c^\infty(U,\mathbb K^r)
\longrightarrow C_c^\infty(U,\mathbb K^r),
\qquad
Pv=\displaystyle\sum_{|\alpha|\leq m}A_\alpha(x)D^\alpha v,
\] be an elliptic system with smooth coefficients. Given $x_0\in U$ and $s\in\mathbb R$, there exist an open subset $V$, with $x_0\in V\Subset U$, and a constant $C_s>0$ such that \[
\|v\|_{H^{s+m}(\mathbb R^n,\mathbb K^r)}
\leq
C_s\bigl(
\|Pv\|_{H^s(\mathbb R^n,\mathbb K^r)}
+
\|v\|_{H^{s+m-1}(\mathbb R^n,\mathbb K^r)}
\bigr)
\] for every $v\in C_c^\infty(V,\mathbb K^r)$, extended by zero outside $V$. \end{proposition}

\begin{proof} We freeze only the principal coefficients at $x_0$ and set \[
P_{x_0}^{(m)}
:=
\displaystyle\sum_{|\alpha|=m}A_\alpha(x_0)D^\alpha.
\] Lemma~\ref{lem:estimacion-eliptica-coeficientes-constantes-orden-arbitrario} provides a constant $C_0>0$ such that \[
\|v\|_{H^{s+m}(\mathbb R^n,\mathbb K^r)}
\leq
C_0\bigl(
\|P_{x_0}^{(m)}v\|_{H^s(\mathbb R^n,\mathbb K^r)}
+
\|v\|_{H^{s+m-1}(\mathbb R^n,\mathbb K^r)}
\bigr).
\] Write \[
P=P_{x_0}^{(m)}+R_1+R_0,
\] where \[
R_1
:=
\displaystyle\sum_{|\alpha|=m}
\bigl(A_\alpha(x)-A_\alpha(x_0)\bigr)D^\alpha
\] and $R_0$ collects the terms of order at most $m-1$. For every $W\Subset U$, the operators $P_{x_0}^{(m)}$, $R_1$, and $R_0$ send $C_c^\infty(W,\mathbb K^r)$ into $C_c^\infty(U,\mathbb K^r)$.

We will apply the refined multiplication estimate, Lemma~\ref{lem:regularidad-intermedia-multiplicacion-afinada}, after extending the coefficients to all of $\mathbb R^n$. Let us make this localization precise. Choose $\rho_0>0$ with $\overline{B(x_0,\rho_0)}\subseteq U$, and define $C_\alpha(\rho_0)$ and $C_*(\rho_0)$ as in the preceding proof. Identities \eqref{eq:variacion-coeficientes-principales-hadamard} and \eqref{eq:cota-coeficiente-principal-radio} show that, for $0<\rho<\rho_0$, \begin{equation}
\sum_{|\alpha|=m}
\sup_{x\in B(x_0,\rho)}
\|A_\alpha(x)-A_\alpha(x_0)\|
\leq \rho C_*(\rho_0).
\label{eq:cota-suma-coeficientes-principales-radio}
\end{equation} Let $C_1>0$ be a constant bounding the product of the matrix constant in Lemma~\ref{lem:regularidad-intermedia-multiplicacion-afinada} and the norms of the operators $D^\alpha\colon H^{s+m}(\mathbb R^n,\mathbb K^r)\longrightarrow
H^s(\mathbb R^n,\mathbb K^r)$, $|\alpha|=m$. This constant depends on $n,s,m$ and the rank $r$, but not on the radius $\rho$, $\theta$, or the coefficients. Choose $\rho\in(0,\rho_0)$ so that $C_0C_1\rho C_*(\rho_0)\leq\displaystyle\frac{1}{2}$, and take $V=B(x_0,\displaystyle\frac{\rho}{2})$. Choose $\theta\in C_c^\infty(B(x_0,\rho))$ with $0\leq\theta\leq1$ and $\theta=1$ on a neighborhood of $\overline V$. If $v\in C_c^\infty(V,\mathbb K^r)$, we may replace each coefficient of $R_1$ by \[
B_\alpha(x)
:=
\theta(x)\bigl(A_\alpha(x)-A_\alpha(x_0)\bigr),
\] interpreted as zero outside $B(x_0,\rho)$. This coefficient is smooth and agrees with $A_\alpha-A_\alpha(x_0)$ on a neighborhood of the support of $v$. Moreover, since $0\leq\theta\leq1$ and its support is contained in $B(x_0,\rho)$, estimate~\eqref{eq:cota-suma-coeficientes-principales-radio} gives \[
\sum_{|\alpha|=m}\|B_\alpha\|_{L^\infty(\mathbb R^n)}
\leq \rho C_*(\rho_0).
\] Lemma~\ref{lem:regularidad-intermedia-multiplicacion-afinada}, applied componentwise, and continuity of $D^\alpha\colon H^{s+m}(\mathbb R^n,\mathbb K^r)\longrightarrow
H^s(\mathbb R^n,\mathbb K^r)$ give \[
\|R_1v\|_{H^s(\mathbb R^n,\mathbb K^r)}
\leq
C_1\rho C_*(\rho_0)\|v\|_{H^{s+m}(\mathbb R^n,\mathbb K^r)}
+C_2\|v\|_{H^{s+m-1}(\mathbb R^n,\mathbb K^r)}.
\] Here the second term collects the remainders that lose one derivative; their constants may depend on derivatives of the coefficients and on $\theta$, but not on $v$.

To handle $R_0$, we extend its coefficients in the same way using $\theta$ and then by zero to $\mathbb R^n$. Since each summand has the form $B_\alpha D^\alpha v$ with $|\alpha|\leq m-1$, Lemma~\ref{lem:regularidad-intermedia-multiplicacion-afinada}, together with Proposition~\ref{prop:regularidad-intermedia-propiedades-escala-hilbertiana}, gives \[
\|R_0v\|_{H^s(\mathbb R^n,\mathbb K^r)}
\leq
C_3\|v\|_{H^{s+m-1}(\mathbb R^n,\mathbb K^r)}.
\]

The identity $P_{x_0}^{(m)}v=Pv-R_1v-R_0v$ implies \[
\|v\|_{H^{s+m}(\mathbb R^n,\mathbb K^r)}
\leq
C_0\|Pv\|_{H^s(\mathbb R^n,\mathbb K^r)}
+C_0C_1\rho C_*(\rho_0)\|v\|_{H^{s+m}(\mathbb R^n,\mathbb K^r)}
+C_4\|v\|_{H^{s+m-1}(\mathbb R^n,\mathbb K^r)}.
\] Our choice of radius satisfies $C_0C_1\rho C_*(\rho_0)\leq\displaystyle\frac{1}{2}$; if $C_*(\rho_0)=0$, any radius for which the preceding localizations are defined suffices. This condition explains quantitatively why the term $R_1$, although it retains order $m$, can be absorbed into the left-hand side. The coefficients of $R_0$ have not been made small: the loss of one order is what allows them to remain in the last term. This gives the asserted estimate. \end{proof}

\begin{proposition}[Higher-order local estimate] \label{prop:estimacion-eliptica-local-orden-superior} Under the hypotheses of Lemma~\ref{lem:estimacion-eliptica-local-sistemas}, for each $k\in\mathbb N_0$ the open subset $V$ can be shrunk, if necessary, so that there exists $C_k>0$ such that \[
\|u\|_{H^{k+m}(\mathbb R^n,\mathbb K^r)}
\leq
C_k\bigl(\|Pu\|_{H^k(\mathbb R^n,\mathbb K^r)}+\|u\|_{L^2(\mathbb R^n,\mathbb K^r)}\bigr)
\] for every $u\in C_c^\infty(V,\mathbb K^r)$, extended by zero outside $V$. \end{proposition}

\begin{proof} Apply Proposition~\ref{prop:estimacion-eliptica-local-orden-arbitrario} with $s=k$. After shrinking $V$ if necessary, we obtain \[
\|u\|_{H^{k+m}(\mathbb R^n,\mathbb K^r)}
\leq
C_1\bigl(
\|Pu\|_{H^k(\mathbb R^n,\mathbb K^r)}
+
\|u\|_{H^{k+m-1}(\mathbb R^n,\mathbb K^r)}
\bigr).
\] Lemma~\ref{lem:interpolacion-sobolev-orden-inferior}, applied with orders $k+m-1$ and $k+m$, provides, for every $\varepsilon>0$, a constant $C_\varepsilon>0$ such that \[
\|u\|_{H^{k+m-1}(\mathbb R^n,\mathbb K^r)}
\leq
\varepsilon\|u\|_{H^{k+m}(\mathbb R^n,\mathbb K^r)}
+C_\varepsilon\|u\|_{L^2(\mathbb R^n,\mathbb K^r)}.
\] Choose $\varepsilon$ so that $C_1\varepsilon<1$ and absorb the corresponding term into the left-hand side. We obtain \[
\|u\|_{H^{k+m}(\mathbb R^n,\mathbb K^r)}
\leq
C_k\bigl(
\|Pu\|_{H^k(\mathbb R^n,\mathbb K^r)}
+
\|u\|_{L^2(\mathbb R^n,\mathbb K^r)}
\bigr),
\] as desired. \end{proof}

\subsection{Globalization on a closed manifold}

The local estimates give a global estimate by means of a finite partition of unity. Commutators with cutoff functions have lower order, so the interpolation inequality controls the terms introduced by patching.

\begin{theorem}[Global elliptic estimate] \label{teo:estimacion-eliptica-global-haces} Let $(M,\mathbf{g})$ be a closed Riemannian manifold, let $\mathbf{E},\mathbf{F}\to M$ be smooth vector bundles with bundle metrics $\mathbf{h}_{\mathbf{E}},\mathbf{h}_{\mathbf{F}}$, and let $P\in\mathbf{PDO}^{(m)}(\mathbf{E},\mathbf{F})$ be elliptic. For every $k\in\mathbb N_0$, there exists $C_k>0$ such that \[
\|\mathbf{u}\|_{H^{k+m}(M,\mathbf{E})}
\leq
C_k\bigl(\|P\mathbf{u}\|_{H^k(M,\mathbf{F})}+\|\mathbf{u}\|_{L^2(M,\mathbf{E})}\bigr)
\] for every $\mathbf{u}\in\Gamma(\mathbf{E})$. Moreover, $P$ extends uniquely to a continuous linear operator \[
P\colon H^{k+m}(M,\mathbf{E})\longrightarrow H^k(M,\mathbf{F}).
\] \end{theorem}

\begin{proof} Continuity of $P\colon H^{k+m}(M,\mathbf{E})\longrightarrow H^k(M,\mathbf{F})$ follows from the local expression in Corollary~\ref{cor:pdo-expresion-con-marcos}. Each derivative of order at most $k$ of $P\mathbf{u}$ is a finite sum of smooth coefficients multiplied by derivatives of the components of $\mathbf{u}$ of order at most $k+m$. Lemma~\ref{lema:meyers-serrin-haz-E}, applied with $p=2$, locally identifies the norm defined using covariant derivatives with the sum of the Sobolev norms of the components. This allows us to sum the estimates over a finite cover.

If $m=0$, the operator $P$ is a smooth bundle homomorphism, and ellipticity means that each $P_x$ is invertible. In local frames, the inverse matrix formula shows that $P^{-1}$ is smooth. The same Leibniz estimate, applied to $P^{-1}$ on a finite cover, gives $\|\mathbf u\|_{H^k(M,\mathbf E)}\leq C_k\|P\mathbf u\|_{H^k(M,\mathbf F)}$. This also proves the reverse inequality at order zero.

Now suppose that $m\geq1$. For the reverse inequality, choose a finite cover $(U_a)_{a=1}^N$ by charts over which $\mathbf{E}$ and $\mathbf{F}$ admit smooth frames, and a subordinate partition of unity $(\chi_a)_{a=1}^N$. Refine the cover so that $\operatorname{supp}(\chi_a)$ lies in an open subset on which Proposition~\ref{prop:estimacion-eliptica-local-orden-superior} applies. In coordinates and components, \[
\|\chi_a\mathbf{u}\|_{H^{k+m}(M,\mathbf{E})}
\leq
C_a\bigl(\|P(\chi_a\mathbf{u})\|_{H^k(M,\mathbf{F})}+\|\chi_a\mathbf{u}\|_{L^2(M,\mathbf{E})}\bigr).
\] Since \[
P(\chi_a\mathbf{u})=\chi_aPu+[P,\chi_a]\mathbf{u}
\] and $[P,\chi_a]$ has order at most $m-1$, \[
\|[P,\chi_a]\mathbf{u}\|_{H^k(M,\mathbf{F})}
\leq
C_a'\|\mathbf{u}\|_{H^{k+m-1}(M,\mathbf{E})}.
\] Summing over $a$ and using Lemma~\ref{lema:meyers-serrin-haz-E} once again, \[
\|\mathbf{u}\|_{H^{k+m}(M,\mathbf{E})}
\leq
C\bigl(\|P\mathbf{u}\|_{H^k(M,\mathbf{F})}+\|\mathbf{u}\|_{H^{k+m-1}(M,\mathbf{E})}+\|\mathbf{u}\|_{L^2(M,\mathbf{E})}\bigr).
\] Specify the passage from local to global interpolation. Put $\ell=k+m$ and let $u_a$ be the coordinate components of $\chi_a\mathbf u$, extended by zero to $\mathbb R^n$. Norm equivalence in the charts and Leibniz's formula give constants $A,B,D>0$, independent of $\mathbf u$, such that
\[
\begin{split}
\|\mathbf u\|_{H^{\ell-1}(M,\mathbf E)}&\leq A\sum_{a=1}^N\|u_a\|_{H^{\ell-1}},\\
\sum_{a=1}^N\|u_a\|_{H^\ell}&\leq B\|\mathbf u\|_{H^\ell(M,\mathbf E)},\\
\sum_{a=1}^N\|u_a\|_{L^2}&\leq D\|\mathbf u\|_{L^2(M,\mathbf E)}.
\end{split}
\]
Component norms are taken on $\mathbb R^n$. Apply Lemma~\ref{lem:interpolacion-sobolev-orden-inferior} to each $u_a$ with parameter $\delta>0$ and sum to obtain
\[
\|\mathbf u\|_{H^{\ell-1}(M,\mathbf E)}
\leq AB\delta\|\mathbf u\|_{H^\ell(M,\mathbf E)}
      +AD C_\delta\|\mathbf u\|_{L^2(M,\mathbf E)}.
\]
Taking $\delta=\displaystyle\frac{\varepsilon}{AB}$ gives the global inequality \[
\|\mathbf{u}\|_{H^{k+m-1}(M,\mathbf{E})}
\leq
\varepsilon\|\mathbf{u}\|_{H^{k+m}(M,\mathbf{E})}+C_\varepsilon\|\mathbf{u}\|_{L^2(M,\mathbf{E})}.
\] Choose $\varepsilon$ so that the first term can be absorbed into the left-hand side. This proves the estimate for smooth sections. Continuity of $P$ and density of $\Gamma(\mathbf{E})$ in $H^{k+m}(M,\mathbf{E})$ yield the indicated extension and extend the inequality to all of $H^{k+m}(M,\mathbf{E})$. \end{proof}

\begin{corollary}[Equivalence with the graph norm] \label{cor:equivalencia-norma-grafica-operador-eliptico} Let $(M,\mathbf{g})$ be a closed Riemannian manifold, let $\mathbf{E},\mathbf{F}\to M$ be smooth vector bundles with bundle metrics $\mathbf{h}_{\mathbf{E}},\mathbf{h}_{\mathbf{F}}$, and let $P\in\mathbf{PDO}^{(m)}(\mathbf{E},\mathbf{F})$ be elliptic. For each $k\in\mathbb N_0$, there exist constants $c_P,C_P>0$, depending on $P$, $k$, and the fixed geometric data, such that \[
c_P\|\mathbf{u}\|_{H^{k+m}(M,\mathbf{E})}
\leq\|\mathbf{u}\|_{L^2(M,\mathbf{E})}+\|P\mathbf{u}\|_{H^k(M,\mathbf{F})}
\leq C_P\|\mathbf{u}\|_{H^{k+m}(M,\mathbf{E})}
\] for every $\mathbf{u}\in H^{k+m}(M,\mathbf{E})$. \end{corollary}

\begin{proof} Theorem~\ref{teo:estimacion-eliptica-global-haces} gives one inequality. The other follows from continuity of $P\colon H^{k+m}(M,\mathbf{E})\longrightarrow H^k(M,\mathbf{F})$ and the continuous inclusion $H^{k+m}(M,\mathbf{E})\hookrightarrow L^2(M,\mathbf{E})$. \end{proof}

In the language of operators on bundles, the preceding estimate agrees with the fundamental elliptic inequality in \cite[Exercise 9.12(b)]{BleeckerBoossIndex}. An alternative construction of elliptic regularity, based on microlocal tools that have not yet been defined, will be developed in Chapter~\ref{cap:operadores-pseudodiferenciales}.

\section{Elliptic regularity}

The operator $P$ acts on distributions through the construction in Chapter~\ref{cap:distribuciones-en-haces-vectoriales}; see Proposition~\ref{prop: pdo induce operador distribuciones}. Thus the equation $P\mathbf{u}=\mathbf{f}$ makes sense for $\mathbf{u}\in\mathcal D'(M,\mathbf{E})$ under the definition of distributional derivative already established.

\begin{lemma}[Local gain of one derivative] \label{lem:ganancia-local-una-derivada-sistema-eliptico} Let $\Omega\subseteq\mathbb R^n$ be open, let \[
P=\displaystyle\sum_{|\alpha|\leq m}A_\alpha(x)D^\alpha
\] be an elliptic system with smooth coefficients, and let $t\in\mathbb Z$. If \[
v\in H^t_{\mathrm{loc}}(\Omega,\mathbb K^r)
\qquad\text{and}\qquad
Pv\in H^{t-m+1}_{\mathrm{loc}}(\Omega,\mathbb K^r),
\] then \[
v\in H^{t+1}_{\mathrm{loc}}(\Omega,\mathbb K^r).
\] \end{lemma}

\begin{proof} Fix $x_0\in\Omega$. Apply Proposition~\ref{prop:estimacion-eliptica-local-orden-arbitrario} with the index \[
s:=t-m+1.
\] We obtain an open subset $V_0$, with $x_0\in V_0\Subset\Omega$, and a constant $C>0$ such that \[
\|z\|_{H^{t+1}(\mathbb R^n,\mathbb K^r)}
\leq
C\bigl(
\|Pz\|_{H^{t-m+1}(\mathbb R^n,\mathbb K^r)}
+
\|z\|_{H^t(\mathbb R^n,\mathbb K^r)}
\bigr)
\] for every $z\in C_c^\infty(V_0,\mathbb K^r)$.

Choose open subsets $W,V$ such that \[
x_0\in W\Subset V\Subset V_0
\] and a function $\chi\in C_c^\infty(V)$ equal to one on a neighborhood of $\overline W$. Define \[
w:=\chi v,
\] extended by zero outside $\Omega$. Since $v\in H^t_{\mathrm{loc}}(\Omega,\mathbb K^r)$, Definition~\ref{def:regularidad-intermedia-sobolev-local-orden-real} gives $w\in H^t(\mathbb R^n,\mathbb K^r)$ and $\operatorname{supp}(w)\subseteq\operatorname{supp}(\chi)\Subset V$.

In the distributional sense, \[
Pw=\chi Pv+[P,\chi]v.
\] Let us justify the regularity of both summands. The hypothesis on $Pv$ gives $\chi Pv\in H^{t-m+1}(\mathbb R^n,\mathbb K^r)$. On the other hand, $[P,\chi]$ is a differential system of order at most $m-1$. Choose $\eta\in C_c^\infty(V_0)$ such that $\eta=1$ on a neighborhood of $\operatorname{supp}(\chi)$. Locality of differential operators allows us to write \[
[P,\chi]v=[P,\chi](\eta v).
\] Since $\eta v\in H^t(\mathbb R^n,\mathbb K^r)$, continuity of differentiation and multiplication by smooth coefficients, Proposition~\ref{prop:regularidad-intermedia-propiedades-escala-hilbertiana} and Lemma~\ref{lem:regularidad-intermedia-multiplicacion-afinada}, gives \[
[P,\chi]v\in H^{t-m+1}(\mathbb R^n,\mathbb K^r).
\] Therefore \[
Pw\in H^{t-m+1}(\mathbb R^n,\mathbb K^r).
\]

Let $J_\varepsilon$ be the mollifier in Proposition~\ref{prop:regularidad-intermedia-molificadores-escala-sobolev} and set $w_\varepsilon:=J_\varepsilon w$. Since $\operatorname{supp}(w)\Subset V$, there exists $\varepsilon_0>0$ such that \[
\operatorname{supp}(w_\varepsilon)\subseteq V_0
\qquad
\text{for }0<\varepsilon<\varepsilon_0.
\] Moreover, $w_\varepsilon\in C_c^\infty(V_0,\mathbb K^r)$, so the preceding elliptic estimate applies: \[
\|w_\varepsilon\|_{H^{t+1}(\mathbb R^n,\mathbb K^r)}
\leq
C\bigl(
\|Pw_\varepsilon\|_{H^{t-m+1}(\mathbb R^n,\mathbb K^r)}
+
\|w_\varepsilon\|_{H^t(\mathbb R^n,\mathbb K^r)}
\bigr).
\]

To control the first term, we extend the coefficients of $P$ outside $V_0$. Choose $\vartheta\in C_c^\infty(\Omega)$ equal to one on a neighborhood of $\overline{V_0}$ and replace $A_\alpha$ with the extension by zero of $\vartheta A_\alpha$ to $\mathbb R^n$. This replacement does not change $Pw$, $Pw_\varepsilon$, or the commutator on their supports. Then \[
Pw_\varepsilon
=
P(J_\varepsilon w)
=
J_\varepsilon(Pw)+[P,J_\varepsilon]w.
\] Proposition~\ref{prop:regularidad-intermedia-molificadores-escala-sobolev} and Corollary~\ref{cor:regularidad-intermedia-conmutador-pdo-molificador} give, uniformly for $0<\varepsilon<\varepsilon_0$, \[
\|Pw_\varepsilon\|_{H^{t-m+1}(\mathbb R^n,\mathbb K^r)}
\leq
\|Pw\|_{H^{t-m+1}(\mathbb R^n,\mathbb K^r)}
+C_1\|w\|_{H^t(\mathbb R^n,\mathbb K^r)}.
\] Also, $\|w_\varepsilon\|_{H^t(\mathbb R^n,\mathbb K^r)}\leq\|w\|_{H^t(\mathbb R^n,\mathbb K^r)}$. Consequently, \[
\|w_\varepsilon\|_{H^{t+1}(\mathbb R^n,\mathbb K^r)}
\leq
C_2\bigl(
\|Pw\|_{H^{t-m+1}(\mathbb R^n,\mathbb K^r)}
+
\|w\|_{H^t(\mathbb R^n,\mathbb K^r)}
\bigr),
\] and the right-hand side is independent of $\varepsilon$.

The space $H^{t+1}(\mathbb R^n,\mathbb K^r)$ is a Hilbert space. Thus there exist a sequence $\varepsilon_j\to0$ and an element $z\in H^{t+1}(\mathbb R^n,\mathbb K^r)$ such that $w_{\varepsilon_j}\rightharpoonup z$ weakly in $H^{t+1}(\mathbb R^n,\mathbb K^r)$. The continuous inclusion $H^{t+1}(\mathbb R^n,\mathbb K^r)\hookrightarrow H^t(\mathbb R^n,\mathbb K^r)$ implies that the same sequence converges weakly to $z$ in $H^t(\mathbb R^n,\mathbb K^r)$. However, Proposition~\ref{prop:regularidad-intermedia-molificadores-escala-sobolev} gives $w_{\varepsilon_j}\to w$ strongly in $H^t(\mathbb R^n,\mathbb K^r)$. Uniqueness of the weak limit shows that $z=w$. Consequently, $w\in H^{t+1}(\mathbb R^n,\mathbb K^r)$.

Since $\chi=1$ on a neighborhood of $\overline W$, we have $v=w$ on $W$. Thus $v\in H^{t+1}_{\mathrm{loc}}(W,\mathbb K^r)$ on a neighborhood of $x_0$. Since $x_0$ was arbitrary, the localization criterion in Lemma~\ref{lem:regularidad-intermedia-criterio-localizacion-sobolev} yields $v\in H^{t+1}_{\mathrm{loc}}(\Omega,\mathbb K^r)$. \end{proof}

\begin{theorem}[Local elliptic regularity] \label{teo:regularidad-eliptica-local-operadores-haces} Let $M$ be a smooth manifold without boundary, let $\mathbf{E},\mathbf{F}\to M$ be smooth vector bundles of finite rank, let $P\in\mathbf{PDO}^{(m)}(\mathbf{E},\mathbf{F})$ be elliptic, and let $k\in\mathbb N_0$. If $\mathbf{u}\in\mathcal D'(M,\mathbf{E})$ and $P\mathbf{u}\in H^k_{\mathrm{loc}}(M,\mathbf{F})$, then $\mathbf{u}\in H^{k+m}_{\mathrm{loc}}(M,\mathbf{E})$. In particular, if $P\mathbf{u}$ is smooth, then $\mathbf{u}$ is smooth. \end{theorem}

\begin{proof} If $m=0$, ellipticity gives a smooth inverse homomorphism $P^{-1}$. Multiplication of distributions by its smooth coefficients is compatible with composition, so $\mathbf u=P^{-1}(P\mathbf u)$. The Leibniz rule in each chart proves the local continuity of $P^{-1}$ on $H^k(M,\mathbf F)$ and its preservation of smoothness. Both conclusions follow directly in this case.

Suppose that $m\geq1$ and fix $x_0\in M$. Choose a regular coordinate ball $(U,\varphi)$ around $x_0$, with $\Omega:=\varphi(U)\subseteq\mathbb R^n$, and smooth frames \[
\mathbf{e}_1,\dots,\mathbf{e}_r
\qquad\text{and}\qquad
\mathbf{f}_1,\dots,\mathbf{f}_r
\] for $\mathbf{E}|_U$ and $\mathbf{F}|_U$, respectively. The ranks agree because the principal symbol of an elliptic operator is an isomorphism between fibers. Shrinking $U$ if necessary, we may assume that $U$ is relatively compact and that $x_0$ corresponds to the point $y_0:=\varphi(x_0)$.

Let $q\in C^\infty(\Omega)$ be the positive coordinate weight of the Riemannian measure, characterized by the change-of-variables formula \[
\int_U \eta\,d\lambda_{\mathbf{g}}
=
\int_\Omega (\eta\circ\varphi^{-1})(y)q(y)\,d\lambda_n(y)
\] for every integrable function $\eta$ supported in $U$. Let $\mathbf{e}^1,\dots,\mathbf{e}^r$ and $\mathbf{f}^1,\dots,\mathbf{f}^r$ be the dual frames of $\mathbf{E}^*|_U$ and $\mathbf{F}^*|_U$, respectively. The local representation of $\mathbf{u}$ is the vector-valued distribution \[
\widehat u=(u^1,\dots,u^r)\in\mathcal D'(\Omega,\mathbb K^r)
\] defined by \[
\langle u^a,\psi\rangle
:=\mathbf{u}\left(((q^{-1}\psi)\circ\varphi)\mathbf{e}^a\right),
\qquad
\psi\in C_c^\infty(\Omega),
\] where the test section is extended by zero outside $U$. Define $\widehat{Pu}=((Pu)^1,\dots,(Pu)^r)$ by \[
 \langle (Pu)^b,\psi\rangle
 :=
 P\mathbf{u}\left(((q^{-1}\psi)\circ\varphi)\mathbf{f}^b\right),
 \qquad
 \psi\in C_c^\infty(\Omega),\quad 1\leq b\leq r,
\] again extending the test section by zero. This normalization by $q^{-1}$ makes the components correspond to Lebesgue measure on $\Omega$.

In the chosen frames, the operator has the form \[
P_U
=
\displaystyle\sum_{|\alpha|\leq m}A_\alpha(y)D^\alpha.
\] The distributional extension of $P$ is compatible with this local representation. Let us verify this assertion. If $\Psi=(\Psi_1,\dots,\Psi_r)\in C_c^\infty(\Omega,(\mathbb K^r)^*)$, identify $q^{-1}\Psi$ with the section of $\mathbf{F}^*|_U$ given by $\displaystyle\sum_{b=1}^r((q^{-1}\Psi_b)\circ\varphi)\mathbf{f}^b$. Let $P_U^{\mathsf t}$ be the bilinear formal transpose of the system with respect to Lebesgue measure; its matrices are transposed, not conjugated. Corollary~\ref{cor:formula-local-transpuesto-formal} gives \[
P'\bigl(q^{-1}\Psi\bigr)
=
q^{-1}P_U^{\mathsf t}\Psi.
\] By the definition of the distributional action, \[
\langle\widehat{Pu},\Psi\rangle
=
(P\mathbf{u})(q^{-1}\Psi)
=
 \mathbf{u}\bigl(P'(q^{-1}\Psi)\bigr)
=
\langle\widehat u,P_U^{\mathsf t}\Psi\rangle
=
\langle P_U\widehat u,\Psi\rangle.
\] Thus \[
\widehat{Pu}=P_U\widehat u
\qquad\text{on }\mathcal D'(\Omega,\mathbb K^r).
\]

By Proposition~\ref{prop:formula-local-simbolo-principal}, the principal symbol of the system is \[
p_m(y,\xi)
=
\displaystyle\sum_{|\alpha|=m}A_\alpha(y)\xi^\alpha.
\] Ellipticity of $P$ means that $p_m(y,\xi)$ is invertible for every $y\in\Omega$ and every $\xi\neq0$. Thus $P_U$ is an elliptic system with smooth coefficients.

The hypothesis $P\mathbf{u}\in H^k_{\mathrm{loc}}(M,\mathbf{F})$ and Lemma~\ref{lema:meyers-serrin-haz-E} imply that the components of $P\mathbf{u}$ in the frame $\mathbf{f}_1,\dots,\mathbf{f}_r$ belong to $H^k_{\mathrm{loc}}(\Omega,\mathbb K^r)$. The preceding normalization merely multiplies test functions by the positive smooth function $q^{-1}$, so the local representation satisfies \[
P_U\widehat u
=
\widehat{Pu}
\in H^k_{\mathrm{loc}}(\Omega,\mathbb K^r).
\]

We now need a starting point for the iteration. Choose $\chi_0\in C_c^\infty(\Omega)$ such that $\chi_0=1$ on a neighborhood of $y_0$. By Proposition~\ref{prop:regularidad-intermedia-distribucion-orden-sobolev-negativo}, there exists $N\in\mathbb N$ such that \[
\chi_0\widehat u
\in
H^{-N}(\mathbb R^n,\mathbb K^r).
\] The localization criterion, Lemma~\ref{lem:regularidad-intermedia-criterio-localizacion-sobolev}, shows that \[
\widehat u
\in
H^{-N}_{\mathrm{loc}}(\Omega_0,\mathbb K^r)
\] on some open neighborhood $\Omega_0\Subset\Omega$ of $y_0$ on which $\chi_0=1$.

Apply Lemma~\ref{lem:ganancia-local-una-derivada-sistema-eliptico} repeatedly. In the first step, take $t=-N$. Since \[
-N-m+1\leq k,
\] the continuous inclusion in Proposition~\ref{prop:regularidad-intermedia-propiedades-escala-hilbertiana} gives \[
P_U\widehat u
\in
H^k_{\mathrm{loc}}(\Omega_0)
\hookrightarrow
H^{-N-m+1}_{\mathrm{loc}}(\Omega_0).
\] The derivative-gain lemma yields \[
\widehat u\in H^{-N+1}_{\mathrm{loc}}(\Omega_0).
\] Suppose inductively that, for some integer $t$ with \[
-N\leq t<k+m,
\] we have obtained $\widehat u\in H^t_{\mathrm{loc}}(\Omega_0)$. Since $t\leq k+m-1$, we have $t-m+1\leq k$, and again \[
P_U\widehat u
\in H^k_{\mathrm{loc}}(\Omega_0)
\hookrightarrow
H^{t-m+1}_{\mathrm{loc}}(\Omega_0).
\] Lemma~\ref{lem:ganancia-local-una-derivada-sistema-eliptico} implies \[
\widehat u\in H^{t+1}_{\mathrm{loc}}(\Omega_0).
\] After finitely many steps, we reach \[
\widehat u
\in
H^{k+m}_{\mathrm{loc}}(\Omega_0,\mathbb K^r).
\]

Lemma~\ref{lema:meyers-serrin-haz-E}, applied in the frame $\mathbf{e}_1,\dots,\mathbf{e}_r$, identifies the Sobolev regularity of a section with that of its local components. Thus \[
\mathbf{u}\in H^{k+m}_{\mathrm{loc}}(M,\mathbf{E})
\] on a neighborhood of $x_0$. Since $x_0$ was arbitrary, the first assertion is proved.

Finally, suppose that $P\mathbf{u}$ is smooth. Then $P\mathbf{u}\in H^k_{\mathrm{loc}}(M,\mathbf{F})$ for every $k\in\mathbb N_0$, and the first part gives \[
\mathbf{u}\in H^{k+m}_{\mathrm{loc}}(M,\mathbf{E})
\qquad
\text{for every }k\in\mathbb N_0.
\] Fix $\ell\in\mathbb N_0$ and choose $k$ so that $k+m>\ell+\displaystyle\frac{n}{2}$. In a chart and a local frame, the Euclidean Sobolev embedding gives \[
H^{k+m}_{\mathrm{loc}}(\mathbb R^n,\mathbb K^r)
\hookrightarrow
C^\ell_{\mathrm{loc}}(\mathbb R^n,\mathbb K^r).
\] Thus $\mathbf{u}\in C^\ell(\mathbf{E})$. Since $\ell$ was arbitrary, we conclude that $\mathbf{u}$ is smooth. \end{proof}

\begin{remark} The geometric part of the proof consists of using a chart and local frames to identify the symbol of the system with the intrinsic symbol of the operator. Once this identification is made, the analytic argument is the same for functions and sections: the matrices $A_\alpha$ couple the components, but uniform invertibility of the symbol on the covector sphere gives the same Fourier estimate. \end{remark}

\subsection{Minimal and maximal realizations on a closed manifold}

Let $P\in\mathbf{PDO}^{(m)}(\mathbf{E},\mathbf{F})$ and suppose that $M$ is closed. Define the initial realization \[
P_0\colon \Gamma(\mathbf{E})\subseteq L^2(M,\mathbf{E})\longrightarrow L^2(M,\mathbf{F}),
\qquad P_0\mathbf{u}:=P\mathbf{u}.
\] Since $\Gamma(\mathbf{E})$ is dense, $P_0$ is densely defined. Its closure, when it exists, is called the minimal realization and is denoted by $P_{\min}$. The maximal realization is defined through the distributional extension of $P$: \[
\mathcal D(P_{\max})
:=
\{\mathbf{u}\in L^2(M,\mathbf{E})\mid P\mathbf{u}\in L^2(M,\mathbf{F})\text{ in the distributional sense}\},
\qquad
P_{\max}\mathbf{u}:=P\mathbf{u}.
\]

\begin{theorem}[Domain of an elliptic operator on a closed manifold] \label{cor:dominio-realizacion-eliptica-compacta} Let $(M,\mathbf{g})$ be a closed Riemannian manifold, let $\mathbf{E},\mathbf{F}\to M$ be smooth vector bundles equipped with bundle metrics, and let $P\in\mathbf{PDO}^{(m)}(\mathbf{E},\mathbf{F})$ be elliptic. Then $P_0$ is closable, \[
P_{\min}=P_{\max},
\qquad
\mathcal D(P_{\min})=\mathcal D(P_{\max})=H^m(M,\mathbf{E}),
\] and there exist constants $c_P,C_P>0$, depending on $P$ and the fixed geometric data, such that \[
c_P\|\mathbf{u}\|_{H^m(M,\mathbf{E})}
\leq\|\mathbf{u}\|_{L^2(M,\mathbf{E})}+\|P_{\max}\mathbf{u}\|_{L^2(M,\mathbf{F})}
\leq C_P\|\mathbf{u}\|_{H^m(M,\mathbf{E})}
\] for every $\mathbf{u}\in H^m(M,\mathbf{E})$. In particular, every closed realization of $P$ lying between the minimal and maximal realizations agrees with them. \end{theorem}

\begin{proof} Let $\mathbf{u}\in\mathcal D(P_{\displaystyle\max})$. Then $\mathbf{u}\in L^2(M,\mathbf{E})$ and $P\mathbf{u}\in L^2(M,\mathbf{F})$ distributionally. Theorem~\ref{teo:regularidad-eliptica-local-operadores-haces} with $k=0$ gives $\mathbf{u}\in H^m_{\mathrm{loc}}(M,\mathbf{E})$. Since $M$ is compact, a finite cover by charts turns this local regularity into $\mathbf{u}\in H^m(M,\mathbf{E})$. Thus \[
\mathcal D(P_{\max})\subseteq H^m(M,\mathbf{E}).
\] The reverse inclusion follows from the continuity $P\colon H^m(M,\mathbf{E})\longrightarrow L^2(M,\mathbf{F})$ in Theorem~\ref{teo:estimacion-eliptica-global-haces}. Therefore $\mathcal D(P_{\displaystyle\max})=H^m(M,\mathbf{E})$. Corollary~\ref{cor:equivalencia-norma-grafica-operador-eliptico} shows that this domain is complete in the graph norm; hence $P_{\displaystyle\max}$ is closed. Since $P_0\subseteq P_{\displaystyle\max}$, the initial operator is closable.

Now take $\mathbf{u}\in H^m(M,\mathbf{E})$. By density, there exists a sequence $(\mathbf{u}_j)\subseteq\Gamma(\mathbf{E})$ such that $\mathbf{u}_j\to \mathbf{u}$ in $H^m(M,\mathbf{E})$. Continuity of $P$ gives $P\mathbf{u}_j\to P\mathbf{u}$ in $L^2(M,\mathbf{F})$. Thus $(\mathbf{u}_j,P\mathbf{u}_j)$ converges to $(\mathbf{u},P\mathbf{u})$ in $L^2(M,\mathbf{E})\times L^2(M,\mathbf{F})$, which shows that $\mathbf{u}\in\mathcal D(P_{\min})$ and $P_{\min}\mathbf{u}=P\mathbf{u}$. We obtain $H^m(M,\mathbf{E})\subseteq\mathcal D(P_{\min})$. Since we always have $P_{\min}\subseteq P_{\displaystyle\max}$, the two realizations agree.

Equivalence of the graph norm and the $H^m(M,\mathbf{E})$ norm is Corollary~\ref{cor:equivalencia-norma-grafica-operador-eliptico} with $k=0$. If $\widetilde P$ is any closed realization of $P$ satisfying \[
P_{\min}\subseteq\widetilde P\subseteq P_{\max},
\] then the equalities already proved imply $P_{\min}=\widetilde P=P_{\displaystyle\max}$. This proves the last assertion. \end{proof}

\begin{corollary}[Closed realization of the connection Laplacian] \label{cor:realizacion-cerrada-laplaciano-conexion} Let $(M,\mathbf{g})$ be a closed Riemannian manifold, and let $(\mathbf{E},\mathbf{h}_{\mathbf{E}},\nabla^{\mathbf{E}})$ be a real vector bundle with a bundle metric or a Hermitian bundle with a compatible connection. Consider \[
 L_{\mathbf{E},0}:=(\nabla^{\mathbf{E}})_h^*\nabla^{\mathbf{E}}
 \colon\Gamma(\mathbf{E})\subset L^2(M,\mathbf{E})\longrightarrow L^2(M,\mathbf{E}).
\] Then \begin{equation}
\label{eq:realizaciones-laplaciano-conexion-cerrada}
 \overline{L_{\mathbf{E},0}}
 =L_{\mathbf{E},\min}=L_{\mathbf{E},\max}=L_{\mathbf{E},0}^*,
 \qquad
 \mathcal D(\overline{L_{\mathbf{E},0}})=H^2(M,\mathbf{E}).
\end{equation} This realization is self-adjoint, nonnegative, and agrees with the Friedrichs realization. In particular, $L_{\mathbf{E},0}$ is essentially self-adjoint. \end{corollary}

\begin{proof} The differential operator $L_{\mathbf{E},0}$ is elliptic of order two and formally self-adjoint. Theorem~\ref{cor:dominio-realizacion-eliptica-compacta} gives $L_{\mathbf{E},\min}=L_{\mathbf{E},\displaystyle\max}$ with domain $H^2(M,\mathbf{E})$. By the definition of the adjoint of a densely defined operator and the formal adjoint identity, \[
 \mathbf{v}\in\mathcal D(L_{\mathbf{E},0}^*)
 \iff L_{\mathbf{E},0}\mathbf{v}\in L^2(M,\mathbf{E})\text{ distributionally},
 \qquad L_{\mathbf{E},0}^*\mathbf{v}=L_{\mathbf{E},0}\mathbf{v};
\] hence $L_{\mathbf{E},0}^*=L_{\mathbf{E},\displaystyle\max}$. Since the minimal realization is the closure of the initial operator, we obtain \eqref{eq:realizaciones-laplaciano-conexion-cerrada}, and equality with its adjoint proves self-adjointness.

For $\mathbf{u}\in\Gamma(\mathbf{E})$, \[
 \langle L_{\mathbf{E},0}\mathbf{u},\mathbf{u}\rangle_{L^2(M,\mathbf{E})}
 =\|\nabla^{\mathbf{E}}\mathbf{u}\|_{L^2(M,T^*M\otimes \mathbf{E})}^2\geq0.
\] Taking the closure preserves nonnegativity. The Friedrichs realization is a self-adjoint extension of $L_{\mathbf{E},0}$ contained in $L_{\mathbf{E},0}^*$; since the minimal and maximal realizations agree, there is no distinct intermediate closed extension, and the Friedrichs realization is the one in \eqref{eq:realizaciones-laplaciano-conexion-cerrada}. The equality $\overline{L_{\mathbf{E},0}}=L_{\mathbf{E},0}^*$ is precisely essential self-adjointness of the initial operator. \end{proof}

\section{Powers of the Bochner Laplacian}

The equivalence in Theorem~\ref{teo:sobolev_cerrada_riemanniana} was obtained by integration by parts and commutation identities for covariant derivatives. Elliptic theory allows us to interpret it as an equivalence of graph norms.

\begin{corollary}[Powers of the Laplacian and Sobolev norms] \label{cor:potencias-laplaciano-bochner-dominios-sobolev} Let $(M,\mathbf{g})$ be a closed Riemannian manifold, and let $\mathbf{E}\to M$ be a smooth vector bundle with a real or Hermitian bundle metric $\mathbf{h}_{\mathbf{E}}$ and a compatible connection $\nabla^{\mathbf{E}}$. Let $L_{\mathbf{E}}=(\nabla^{\mathbf{E}})_h^*\nabla^{\mathbf{E}}$ be the realization in \eqref{eq:realizaciones-laplaciano-conexion-cerrada}. For every $q\in\mathbb N$, \[
\mathcal D(L_{\mathbf{E}}^q)=H^{2q}(M,\mathbf{E}),
\] and there exist constants $c_q,C_q>0$ such that \[
c_q\|\mathbf{u}\|_{H^{2q}(M,\mathbf{E})}
\leq
\|\mathbf{u}\|_{L^2(M,\mathbf{E})}+\|L_{\mathbf{E}}^q\mathbf{u}\|_{L^2(M,\mathbf{E})}
\leq
C_q\|\mathbf{u}\|_{H^{2q}(M,\mathbf{E})}.
\] Moreover, \[
\bigcap_{q=1}^\infty\mathcal D(L_{\mathbf{E}}^q)=\Gamma(\mathbf{E}).
\] \end{corollary}

\begin{proof}
Put $P_q=L_{\mathbf E,0}^q$. Proposition~\ref{prop:simbolo-composicion-pdo} and Example~\ref{ej:simbolo-laplaciano-bochner} show that $P_q$ is elliptic of order $2q$. By Theorem~\ref{cor:dominio-realizacion-eliptica-compacta}, its minimal and maximal realizations agree, have domain $H^{2q}(M,\mathbf E)$, and have graph norm equivalent to that space's norm.

The operator $L_{\mathbf E}$ is nonnegative and self-adjoint. The spectral Theorem~\ref{teo:calculo-espectral-autoadjunto-no-negativo} implies that $L_{\mathbf E}^q$ is closed. Smooth sections belong to the domain of every power because $L_{\mathbf E}$ agrees there with $L_{\mathbf E,0}$, which preserves smooth sections. Thus $P_q\subseteq L_{\mathbf E}^q$ and, by closedness, $(P_q)_{\min}\subseteq L_{\mathbf E}^q$.

For the reverse inclusion, let $\mathbf u\in\mathcal D(L_{\mathbf E}^q)$. The inequality $\lambda^{2j}\leq1+\lambda^{2q}$, for $0\leq j\leq q$, shows by spectral calculus that all intermediate powers are defined. Moreover,
\[
L_{\mathbf E}^j\mathbf u\in\mathcal D(L_{\mathbf E}),\qquad 0\leq j<q,
\]
since this domain condition is integrability of the moment of order $2(j+1)$. For each $\mathbf v\in\mathcal D(L_{\mathbf E})$, we have $L_{\mathbf E,0}\mathbf v=L_{\mathbf E}\mathbf v$ distributionally, by the identification of $L_{\mathbf E}$ with the maximal realization. Applying this identity successively to $\mathbf u,L_{\mathbf E}\mathbf u,\ldots,L_{\mathbf E}^{q-1}\mathbf u$ gives
\[
L_{\mathbf E,0}^j\mathbf u=L_{\mathbf E}^j\mathbf u
\quad\text{distributionally},\qquad 0\leq j\leq q.
\]
In particular, $P_q\mathbf u\in L^2(M,\mathbf E)$ and $L_{\mathbf E}^q\subseteq(P_q)_{\max}$. Therefore
\[
(P_q)_{\min}\subseteq L_{\mathbf E}^q\subseteq(P_q)_{\max}=(P_q)_{\min}.
\]
The realizations agree. The domain identity and graph norm equivalence give constants $c_q,C_q>0$ such that
\[
c_q\|\mathbf u\|_{H^{2q}(M,\mathbf E)}
\leq\|\mathbf u\|_{L^2(M,\mathbf E)}+\|L_{\mathbf E}^q\mathbf u\|_{L^2(M,\mathbf E)}
\leq C_q\|\mathbf u\|_{H^{2q}(M,\mathbf E)}.
\]
If $\mathbf u$ belongs to every domain, it belongs to $H^{2q}(M,\mathbf E)$ for every $q$. Local Sobolev embeddings give $C^k$ regularity for every $k$, so $\mathbf u$ is smooth. The converse inclusion was already established.
\end{proof}

\begin{remark} On a complete, noncompact manifold without boundary, local elliptic regularity remains valid, but a global estimate with a single constant requires uniform control of the geometry and the coefficients. Analogous estimates hold under bounded geometry hypotheses; we will return to this point in the corresponding chapter. \end{remark}

The principal symbol constructed here is also the symbol that enters index theory. If $\pi\colon T^*M\longrightarrow M$ is the projection, an elliptic operator determines an isomorphism over $T^*M\setminus0_M$ \[
\pi^*\mathbf{E}\longrightarrow\pi^*\mathbf{F},
\qquad
(x,\xi,v)\longmapsto(x,\xi,\boldsymbol{\sigma}_m(P)(x,\xi)v).
\] The restriction of this isomorphism to the cosphere bundle contains the topological information that will be used later. This covariant formulation agrees with the one developed in \cite[Definition 6.36 and Chapter 9]{BleeckerBoossIndex}.

\chapter{Fractional Sobolev spaces on vector bundles and the heat kernel: an intrinsic approach} \label{cap:sobolev-fraccionario-haces-nucleo-calor}

The integer-order Sobolev spaces studied in the preceding chapters are constructed using the covariant derivatives $\nabla^j\mathbf{u}$. When the order is no longer an integer, no finite family of derivatives by itself gives an intrinsic definition. One alternative is to choose a nonnegative elliptic operator, construct the semigroup it generates, and define its negative powers through a Laplace integral. This is the viewpoint of Bessel potentials.

The semigroup arises in solving the evolution problem \[
\partial_t\mathbf{u}+L_{\mathbf{E}}\mathbf{u}=0,
\qquad
\mathbf{u}(0)=\mathbf{u}_0.
\] Formally, its solution is $\mathbf{u}(t)=e^{-tL_{\mathbf{E}}}\mathbf{u}_0$. Integrating these operators with respect to $t$ with the appropriate weight produces negative powers of $I+L_{\mathbf{E}}$ and, with them, a scale of fractional regularity. Heat thus links a parabolic PDE to the functional calculus of an elliptic operator.

The development proceeds on two levels. In $L^2(M,\mathbf{E})$, we have the Hilbert space structure, closed sesquilinear forms, and spectral calculus for self-adjoint operators. Once the heat semigroup has been constructed on this space, the Kato--Simon inequality allows us to compare it with the scalar semigroup and extend it consistently to $L^p(M,\mathbf{E})$. Bessel potentials are then defined by Bochner integrals. For scalar functions on complete manifolds, this program was developed by Strichartz; see \cite[Theorems 3.5 and 4.2]{LaplacianStr}. Here we formulate it directly for sections and prove the covariant steps involved in the construction.

Throughout this chapter, $(M,\mathbf{g})$ will be a connected, complete Riemannian manifold without boundary. The bundle $\mathbf{E}\to M$ will be a smooth real or complex vector bundle of finite rank, equipped with a bundle metric $\mathbf{h}_{\mathbf{E}}$, Hermitian in the complex case, and $\nabla^{\mathbf{E}}$ will be a connection compatible with $\mathbf{h}_{\mathbf{E}}$. In the real case, spectral theory is applied to the complexification of $\mathbf{E}$; the Laplacian, semigroup, and Bessel potentials of real order preserve real sections. We write $L^p(M,\mathbf{E})$ for the spaces defined with respect to the Riemann--Lebesgue measure $\lambda_{\mathbf{g}}$. On Sobolev classes, $\nabla^j\mathbf{u}$ denotes the weak covariant derivative $\nabla_w^j\mathbf{u}$.

\begin{semblanzaHistorica}{Heat as a tool for regularity} Since Fourier, the heat equation has served to separate and damp frequencies. On a manifold, the semigroup generated by a Laplacian retains this role, even though global translations are no longer available. Spectral theory constructs heat from the operator; its estimates reveal the geometry of the manifold; and integrals of the semigroup produce Bessel potentials and fractional orders. The same object thus connects evolution, regularization, and functional analysis. \end{semblanzaHistorica}

\section{The connection Laplacian on \texorpdfstring{$L^2(M,\mathbf{E})$}{L2(E)}}

We use the nonnegative realization of the Bochner Laplacian from
Definition~\ref{def:operadores-diferenciales-en-haces-laplaciano-de-bochner}.
Covariant derivatives and contractions follow the conventions of
Chapter~\ref{cap:sobolev-haces}, in particular
Theorem~\ref{teo:green-m-esima-covariante-E}. In the scalar case,
$L_{\mathbf E}$ is the realization of $-\Delta_{\mathbf g}$, and
$\partial_tu+L_{\mathbf E}u=0$ is equivalent, for smooth solutions, to
$\partial_tu=\Delta_{\mathbf g}u$.

\subsection{The realization determined by the energy form}

We write \[
L_{\mathbf{E},0}:=\Delta_B\restriction_{\Gamma_c(\mathbf{E})}
=(\nabla^{\mathbf{E}})_h^*\nabla^{\mathbf{E}}
=-\operatorname{tr}_{\mathbf{g}}(\nabla^2),
\qquad
\mathcal D(L_{\mathbf{E},0})=\Gamma_c(\mathbf{E}),
\] for its initial restriction to compactly supported smooth sections. Its distributional action requires no new definition: it is obtained by applying Proposition~\ref{prop: pdo induce operador distribuciones} to the differential operator $L_{\mathbf{E},0}$. In particular, if $\mathbf{u}\in\mathcal D'(M,\mathbf{E})$, the expression $L_{\mathbf{E},0}\mathbf{u}$ has the meaning established in the chapter on distributions on vector bundles.

The local formula for the Bochner Laplacian, proved in the chapter on partial differential operators, gives \[
L_{\mathbf{E},0}\mathbf{u}
=
-\displaystyle\sum_{i=1}^n
\left(
\nabla^{\mathbf{E}}_{\mathbf{e}_i}\nabla^{\mathbf{E}}_{\mathbf{e}_i}\mathbf{u}
-
\nabla^{\mathbf{E}}_{\nabla^{TM}_{\mathbf{e}_i}\mathbf{e}_i}\mathbf{u}
\right)
\] in any local orthonormal frame $(\mathbf{e}_1,\dots,\mathbf{e}_n)$ of $TM$. Example~\ref{ej:simbolo-laplaciano-bochner} shows that, with the differential-symbol convention of Definition~\ref{def:simbolo-principal-pdo-haces}, its principal symbol is $-\lvert\xi\rvert_{\mathbf{g}}^2\operatorname{Id}_{\mathbf{E}_x}$; it is therefore elliptic.

\begin{notation}[Connection Laplacian] \label{def:laplaciano-conexion-capitulo-fraccionario} \index{connection Laplacian@connection Laplacian} In what follows, we call the preceding operator $L_{\mathbf{E},0}$ the \textbf{initial connection Laplacian}. The letter $L_{\mathbf{E}}$ is reserved for its self-adjoint realization on $L^2(M,\mathbf{E})$. \end{notation}

\begin{corollary}[The closed case] \label{cor:caso-cerrado-laplaciano-conexion-capitulo-fraccionario} If, in addition, $M$ is closed, the elliptic results already proved give \[
 \overline{L_{\mathbf{E},0}}=L_{\mathbf{E},\min}=L_{\mathbf{E},\max}=L_{\mathbf{E},0}^*,
 \qquad
 \mathcal D(\overline{L_{\mathbf{E},0}})=H^2(M,\mathbf{E}).
\] The common realization is self-adjoint, nonnegative, and agrees with the Friedrichs realization. In particular, the initial operator $L_{\mathbf{E},0}$ is essentially self-adjoint. \end{corollary}

\begin{proof} This is Corollary~\ref{cor:realizacion-cerrada-laplaciano-conexion}, obtained from elliptic regularity, equality of the minimal and maximal realizations, and formal adjunction. \end{proof}

\begin{proposition} \label{prop:laplaciano-conexion-simetrico-no-negativo} Let $(M,\mathbf{g})$ be a connected, complete Riemannian manifold without boundary, and let $\mathbf{E}\to M$ be a smooth vector bundle of finite rank with a bundle metric and a compatible connection. The operator $L_{\mathbf{E},0}$ is densely defined, symmetric, and nonnegative. For any $\mathbf{u},\mathbf{v}\in\Gamma_c(\mathbf{E})$, we have \[
\langle L_{\mathbf{E},0}\mathbf{u},\mathbf{v}\rangle_{L^2(M,\mathbf{E})}
=
\int_M\langle\nabla^{\mathbf{E}}\mathbf{u},\nabla^{\mathbf{E}}\mathbf{v}\rangle_{\mathbf{g},\mathbf{h}_{\mathbf{E}}}\,d\lambda_{\mathbf{g}}
=
\langle \mathbf{u},L_{\mathbf{E},0}\mathbf{v}\rangle_{L^2(M,\mathbf{E})}.
\] In particular, \[
\langle L_{\mathbf{E},0}\mathbf{u},\mathbf{u}\rangle_{L^2(M,\mathbf{E})}
=
\|\nabla^{\mathbf{E}}\mathbf{u}\|_{L^2(M,T^*M\otimes \mathbf{E})}^2.
\] \end{proposition}

\begin{proof} Density of $\Gamma_c(\mathbf{E})$ in $L^2(M,\mathbf{E})$ was proved in the chapter on Sobolev spaces on bundles. The first equality is the integration-by-parts formula characterizing $(\nabla^{\mathbf{E}})_h^*$; the second follows by interchanging $\mathbf{u}$ and $\mathbf{v}$ and using the Hermitian symmetry of the inner product. Taking $\mathbf{v}=\mathbf{u}$ gives the last identity and hence nonnegativity. \end{proof}

Define on $\Gamma_c(\mathbf{E})$ the form \[
\mathfrak q_{\mathbf{E}}(\mathbf{u},\mathbf{v})
:=
\int_M\langle\nabla^{\mathbf{E}}\mathbf{u},\nabla^{\mathbf{E}}\mathbf{v}\rangle_{\mathbf{g},\mathbf{h}_{\mathbf{E}}}\,d\lambda_{\mathbf{g}}.
\] Its closure has domain $W_0^{1,2}(M,\mathbf{E})$. To identify this space with $W^{1,2}(M,\mathbf{E})$, we need cutoff functions whose gradients tend uniformly to zero. Before constructing them, we give a detailed proof of the existence of a smooth exhaustion function with bounded gradient.

\begin{lemma}[Smoothing Lipschitz functions] \label{lem:suavizacion-global-lipschitz-variedad} Let $(M,\mathbf{g})$ be a connected, complete Riemannian manifold without boundary. Let $f\colon M\longrightarrow\mathbb R$ be an $L$-Lipschitz function with respect to the Riemannian distance, and let $\varepsilon>0$. There exists $h\in C^\infty(M)$ such that \[
|h-f|<\varepsilon
\qquad\text{and}\qquad
|dh|_{\mathbf{g}}\leq L+\varepsilon
\] throughout $M$. \end{lemma}

\begin{proof} If $L=0$, connectedness of $M$ implies that $f$ is constant, and it suffices to take $h=f$. Suppose that $L>0$ and fix \[
\theta:=\min\left\{\frac{1}{4},\frac{\varepsilon}{16L}\right\}.
\] For each $p\in M$, Proposition~\ref{prop:bolas-coordenadas-normales-control-euclidiano}, applied with this value of $\theta$, provides a regular normal coordinate ball $U_p$ and a normal chart $\kappa_p\colon U_p\longrightarrow B_{\mathrm{euc}}(0,r_p)$ such that, for every $z\in B_{\mathrm{euc}}(0,r_p)$ and $w\in\mathbb R^n$, \[
(1-\theta)\|w\|
\leq
\left|d(\kappa_p^{-1})_z w\right|_{\mathbf{g}}
\leq
(1+\theta)\|w\|.
\] These balls cover $M$. Consider the family $\mathcal B$ of regular coordinate balls $V$ for which $\overline V\subseteq U_p$ for some $p\in M$. This family is a basis: if $q\in O$ and $O$ is open, choose $p$ with $q\in U_p$ and, inside the Euclidean open subset $\kappa_p(O\cap U_p)$, two concentric balls centered at $\kappa_p(q)$ whose radii satisfy $0<s<s'$ and whose closed ball of radius $s'$ is contained in $\kappa_p(O\cap U_p)$. The preimage of the ball of radius $s$ is a regular coordinate ball $V$ such that $q\in V$ and $\overline V\subseteq O\cap U_p$, in accordance with Definition~\ref{def:nociones-fundamentales-norma-variedad-suave-label-bola}. Theorem~\ref{teo: variedad es paracompacta} allows us to choose from $\mathcal B$ a countable, locally finite refinement $(V_j)_{j\in\mathbb N}$ of the cover $(U_p)_{p\in M}$. For each $j$, fix a normal ball $U_j$ from the original cover such that $\overline V_j\subseteq U_j$, and denote its normal chart by $\kappa_j\colon U_j\longrightarrow B_{\mathrm{euc}}(0,r_j)$. Let $(\varphi_j)_{j\in\mathbb N}$ be a smooth partition of unity subordinate to $(V_j)_{j\in\mathbb N}$, whose existence is guaranteed by Theorem~\ref{particionesdelaunidad}. Since $\operatorname{supp}(\varphi_j)$ is closed in $M$ and contained in the compact set $\overline V_j$, the numbers \[
C_j:=\|d\varphi_j\|_{L^\infty(M,T^*M)},
\qquad
\delta_j:=\min\left\{\frac{\varepsilon}{4},
\frac{\varepsilon}{2^{j+4}(1+C_j)}\right\}
\] are well defined, with $C_j\geq0$ and $\delta_j>0$.

Fix $j\in\mathbb N$ and write $F_j:=f\circ\kappa_j^{-1}$. If $z,z'\in B_{\mathrm{euc}}(0,r_j)$, the segment joining them remains in this Euclidean ball. The curve $\gamma(t):=\kappa_j^{-1}((1-t)z+tz')$ lies in $U_j$ and satisfies \[
d_{\mathbf{g}}\bigl(\kappa_j^{-1}(z),\kappa_j^{-1}(z')\bigr)
\leq
L_{\mathbf{g}}(\gamma)
\leq
(1+\theta)\|z-z'\|.
\] Since $f$ is $L$-Lipschitz with respect to $d_{\mathbf{g}}$, it follows that \[
|F_j(z)-F_j(z')|
\leq
L(1+\theta)\|z-z'\|.
\] Thus $F_j$ is $L_j$-Lipschitz, where $L_j:=L(1+\theta)$. McShane's extension theorem, Theorem~\ref{teo:extension-mcshane}, gives an extension $\widetilde F_j\colon \mathbb R^n\longrightarrow\mathbb R$ with the same constant $L_j$. Take the mollifier $\rho$ from Theorem~\ref{teo:molificadores-euclidianos} and set \[
\tau_j:=\min\left\{1,
\frac{\delta_j}{2L_jm_1(\rho)}\right\},
\qquad
f_j:=\bigl(\widetilde F_j*\rho_{\tau_j}\bigr)\circ\kappa_j
\quad\text{in }U_j.
\] The estimates in that theorem imply $|f_j-f|\leq\displaystyle\frac{\delta_j}{2}$ on $U_j$ and $\|D(\widetilde F_j*\rho_{\tau_j})\|\leq L_j$. If $q\in U_j$, $z=\kappa_j(q)$, and $v\in T_qM$ satisfies $|v|_{\mathbf{g}}=1$, the first metric inequality gives $\|d\kappa_j(v)\|\leq\displaystyle\frac{1}{1-\theta}$. Therefore \[
|df_j|_{\mathbf{g}}
\leq
\sup_{\substack{v\in T_qM\\|v|_{\mathbf{g}}=1}}
\left|D(\widetilde F_j*\rho_{\tau_j})_z\bigl(d\kappa_j(v)\bigr)\right|
\leq
\frac{L_j}{1-\theta}
=
L\frac{1+\theta}{1-\theta}.
\] Moreover, \[
L\frac{1+\theta}{1-\theta}-L
=
\frac{2L\theta}{1-\theta}
\leq
\frac{\varepsilon}{6}.
\] Indeed, if $\theta=\displaystyle\frac{\varepsilon}{16L}$, use $\theta\leq\displaystyle\frac{1}{4}$; if $\theta=\displaystyle\frac{1}{4}$, then $\varepsilon\geq4L$ and the same bound follows.

Each product $\varphi_jf_j$, defined to be zero outside $U_j$, is smooth because $\operatorname{supp}(\varphi_j)\subseteq V_j\Subset U_j$. Define $h:=\displaystyle\sum_{j=1}^{\infty}\varphi_jf_j$. The sum is locally finite, so $h\in C^\infty(M)$. Since $\displaystyle\sum_{j=1}^{\infty}\varphi_j=1$ and $\varphi_j\geq0$, \[
|h-f|
\leq
\displaystyle\sum_{j=1}^{\infty}\varphi_j|f_j-f|
\leq
\frac{1}{2}\displaystyle\sum_{j=1}^{\infty}\varphi_j\delta_j
\leq
\frac{\varepsilon}{8}
<\varepsilon.
\] The identity $\displaystyle\sum_{j=1}^{\infty}d\varphi_j=0$ allows us to write, without differentiating the Lipschitz function $f$, \[
dh
=
\displaystyle\sum_{j=1}^{\infty}\varphi_jdf_j
+
\displaystyle\sum_{j=1}^{\infty}(f_j-f)d\varphi_j.
\] By the preceding bounds and local finiteness, \[
|dh|_{\mathbf{g}}
\leq
L+\frac{\varepsilon}{6}
+
\frac{1}{2}\displaystyle\sum_{j=1}^{\infty}\delta_jC_j
\leq
L+\frac{\varepsilon}{6}
+
\frac{\varepsilon}{2}\displaystyle\sum_{j=1}^{\infty}2^{-j-4}
<
L+\varepsilon.
\] This proves both assertions. \end{proof}

\begin{remark} The preceding proof does not require the balls $U_j$ to be strongly geodesically convex. The convexity used is that of $\kappa_j(U_j)=B_{\mathrm{euc}}(0,r_j)$: the Euclidean segment determines a curve contained in $U_j$, and its length gives an upper bound for the Riemannian distance. Thus the local control in Proposition~\ref{prop:bolas-coordenadas-normales-control-euclidiano} is exactly the geometric result needed. \end{remark}

\begin{proposition}[Exhaustion function with bounded gradient] \label{prop:agotamiento-suave-gradiente-acotado} Let $(M,\mathbf{g})$ be a connected, complete Riemannian manifold without boundary. There exists a strictly positive proper function $\rho\in C^\infty(M)$ such that \[
|d\rho|_{\mathbf{g}}\leq2
\] throughout $M$. \end{proposition}

\begin{proof} Fix $o\in M$ and consider $r(x):=d_{\mathbf{g}}(o,x)$. The triangle inequality shows that $r$ is $1$-Lipschitz. Since $M$ is complete, the Hopf--Rinow theorem, Theorem~\ref{teo:variedades-riemannianas-hopf-rinow}, implies that the closed balls $\overline{B}_{\mathbf{g}}(o,R)$ are compact; hence $r$ is proper.

Apply Lemma~\ref{lem:suavizacion-global-lipschitz-variedad} to $r$ with $\varepsilon=1$. We obtain $h\in C^\infty(M)$ such that $|h-r|<1$ and $|dh|_{\mathbf{g}}\leq2$. In particular, $|h(o)|<1$ because $r(o)=0$. Set $\rho:=h-h(o)+2$. This function is strictly positive because \[
\rho(x)
\geq
r(x)-1-h(o)+2
\geq
1-h(o)>0,
\], and it is proper: if $\rho(x)\leq R$, then $r(x)\leq R+h(o)-1$, so the sublevel set $\{\rho\leq R\}$ is contained in a closed ball and is closed; by Theorem~\ref{teo:variedades-riemannianas-hopf-rinow}, it is compact. Finally, $d\rho=dh$, so $|d\rho|_{\mathbf{g}}\leq2$. \end{proof}

\begin{lemma}[Cutoff functions on complete manifolds] \label{lem:funciones-corte-completitud-capitulo-calor} Let $(M,\mathbf{g})$ be a connected, complete Riemannian manifold without boundary. \index{cutoff function@cutoff function!on a complete manifold} There exists a sequence $(\chi_j)_{j\in\mathbb N}\subseteq C_c^\infty(M)$ such that $0\leq\chi_j\leq1$, $\chi_j\to1$ uniformly on every compact set, and \[
\|d\chi_j\|_{L^\infty(M,T^*M)}\longrightarrow0.
\] \end{lemma}

\begin{proof} Let $\rho$ be the function in Proposition~\ref{prop:agotamiento-suave-gradiente-acotado}. Choose a nonincreasing function $\eta\in C^\infty(\mathbb R)$ such that $0\leq\eta\leq1$, $\eta(s)=1$ for $s\leq1$, and $\eta(s)=0$ for $s\geq2$, and define \[
\chi_j(x):=\eta\left(\frac{\rho(x)}{j}\right).
\] The support of $\chi_j$ is contained in the sublevel set $\{\rho\leq2j\}$, which is compact because $\rho$ is proper. Since $\rho>0$ and $\eta$ is nonincreasing, the sequence $(\chi_j)$ is increasing. If $K\subseteq M$ is compact and $\displaystyle R_K:=\displaystyle\max_{x\in K}\rho(x)$, then $\chi_j=1$ on $K$ for every integer $\displaystyle j\geq\displaystyle\max\{1,\lceil R_K\rceil\}$. In particular, $\chi_j$ converges pointwise to $1$. Finally, \[
|d\chi_j|_{\mathbf{g}}
=
\frac{1}{j}\left|\eta'\left(\frac{\rho}{j}\right)\right||d\rho|_{\mathbf{g}}
\leq
\frac{2\|\eta'\|_{L^\infty(\mathbb R)}}{j},
\] which proves the asserted convergence. \end{proof}

\begin{proposition} \label{prop:W12-cero-igual-W12-completa} Let $(M,\mathbf{g})$ be a connected, complete Riemannian manifold without boundary, and let $\mathbf{E}\to M$ be a smooth vector bundle of finite rank with a bundle metric and a compatible connection. Then $W_0^{1,2}(M,\mathbf{E})=W^{1,2}(M,\mathbf{E})$. \end{proposition}

\begin{proof} Let $\mathbf{u}\in W^{1,2}(M,\mathbf{E})$, and let $(\chi_j)$ be the functions in the preceding lemma. The weak Leibniz rule gives \[
\nabla^{\mathbf{E}}(\chi_j\mathbf{u})=d\chi_j\otimes \mathbf{u}+\chi_j\nabla^{\mathbf{E}}\mathbf{u}.
\] Since $0\leq\chi_j\leq1$ and $\chi_j\to1$ pointwise, we have \[
|(1-\chi_j)\mathbf{u}|_{\mathbf{h}_{\mathbf{E}}}^2\leq|\mathbf{u}|_{\mathbf{h}_{\mathbf{E}}}^2,
\qquad
|(1-\chi_j)\nabla^{\mathbf{E}}\mathbf{u}|_{\mathbf{g},\mathbf{h}_{\mathbf{E}}}^2\leq|\nabla^{\mathbf{E}}\mathbf{u}|_{\mathbf{g},\mathbf{h}_{\mathbf{E}}}^2.
\] Both dominating functions are integrable. Lebesgue's dominated convergence theorem~\ref{convergencia dominada} implies \[
(1-\chi_j)\mathbf{u}\longrightarrow0\quad\text{in }L^2(M,\mathbf{E})
\] and \[
(1-\chi_j)\nabla^{\mathbf{E}}\mathbf{u}\longrightarrow0
\quad\text{in }L^2(M,T^*M\otimes \mathbf{E}).
\] Moreover, \[
\|d\chi_j\otimes \mathbf{u}\|_{L^2(M,T^*M\otimes \mathbf{E})}
\leq
\|d\chi_j\|_{L^\infty(M,T^*M)}\|\mathbf{u}\|_{L^2(M,\mathbf{E})}
\longrightarrow0.
\] Therefore $\chi_j\mathbf{u}\to \mathbf{u}$ in $W^{1,2}(M,\mathbf{E})$. For each $j$, the Meyers--Serrin theorem for bundles, Theorem~\ref{meyers-serrin-haz}, approximates $\chi_j\mathbf{u}$ in $W^{1,2}(M,\mathbf{E})$ by smooth sections $\mathbf{u}_{j,k}$. Choose $\zeta_j\in C_c^\infty(M)$ equal to $1$ on a neighborhood of $\operatorname{supp}(\chi_j\mathbf{u})$. Since multiplication by $\zeta_j$ is continuous on $W^{1,2}(M,\mathbf{E})$, we have $\zeta_j\mathbf{u}_{j,k}\to\chi_j\mathbf{u}$, and now $\zeta_j\mathbf{u}_{j,k}\in\Gamma_c(\mathbf{E})$. Thus $\mathbf{u}\in W_0^{1,2}(M,\mathbf{E})$. \end{proof}

The closure of $\mathfrak q_{\mathbf{E}}$ is therefore a densely defined, symmetric, nonnegative, closed sesquilinear form with domain $W^{1,2}(M,\mathbf{E})$.

\begin{theorem}[Friedrichs realization of the connection Laplacian] \label{teo:realizacion-friedrichs-laplaciano-conexion} Let $(M,\mathbf{g})$ be a connected, complete Riemannian manifold without boundary, and let $\mathbf{E}\to M$ be a smooth vector bundle of finite rank with a bundle metric and a compatible connection. There exists a unique nonnegative self-adjoint operator \[
L_{\mathbf{E}}\colon \mathcal D(L_{\mathbf{E}})\subseteq L^2(M,\mathbf{E})\longrightarrow L^2(M,\mathbf{E})
\] such that $\mathcal D(L_{\mathbf{E}})\subseteq W^{1,2}(M,\mathbf{E})$ and \[
\mathfrak q_{\mathbf{E}}(\mathbf{u},\mathbf{v})=\langle L_{\mathbf{E}}\mathbf{u},\mathbf{v}\rangle_{L^2(M,\mathbf{E})},
\qquad \mathbf{u}\in\mathcal D(L_{\mathbf{E}}),\quad \mathbf{v}\in W^{1,2}(M,\mathbf{E}).
\] Its domain is \[
\mathcal D(L_{\mathbf{E}})
=
\left\{\mathbf{u}\in W^{1,2}(M,\mathbf{E})\middle|
\begin{array}{l}
\text{there exists }\mathbf{f}\in L^2(M,\mathbf{E})\text{ such that}\\[2pt]
\mathfrak q_{\mathbf{E}}(\mathbf{u},\mathbf{v})=\langle \mathbf{f},\mathbf{v}\rangle_{L^2(M,\mathbf{E})}
\text{ for every }\mathbf{v}\in W^{1,2}(M,\mathbf{E})
\end{array}
\right\},
\] and, in this case, $L_{\mathbf{E}}\mathbf{u}=\mathbf{f}$. Moreover, $L_{\mathbf{E}}$ extends $L_{\mathbf{E},0}$. \end{theorem}

\begin{proof} Existence, uniqueness, and the characterization of the domain follow from Theorem~\ref{teo:representacion-formas-cerradas}. If $\mathbf{u}\in\Gamma_c(\mathbf{E})$, Proposition~\ref{prop:laplaciano-conexion-simetrico-no-negativo} gives \[
\mathfrak q_{\mathbf{E}}(\mathbf{u},\mathbf{v})=\langle L_{\mathbf{E},0}\mathbf{u},\mathbf{v}\rangle_{L^2(M,\mathbf{E})}
\] for $\mathbf{v}\in\Gamma_c(\mathbf{E})$. Both sides depend continuously on $\mathbf{v}$ in the $W^{1,2}(M,\mathbf{E})$ norm, and $\Gamma_c(\mathbf{E})$ is dense in that space by Proposition~\ref{prop:W12-cero-igual-W12-completa}. The identity extends to all of $\mathbf{v}\in W^{1,2}(M,\mathbf{E})$, which proves that $\mathbf{u}\in\mathcal D(L_{\mathbf{E}})$ and $L_{\mathbf{E}}\mathbf{u}=L_{\mathbf{E},0}\mathbf{u}$. \end{proof}

\begin{theorem}[Essential self-adjointness] \label{teo:sobolev-fraccionario-y-nucleo-calor-variedad-riemanniana-completa-esencialmente-auto} Let $(M,\mathbf{g})$ be a connected, complete Riemannian manifold without boundary, and let $\mathbf{E}\to M$ be a smooth vector bundle of finite rank with a bundle metric and a compatible connection. \index{connection Laplacian@connection Laplacian!essential self-adjointness} The operator $L_{\mathbf{E},0}$ is essentially self-adjoint on $L^2(M,\mathbf{E})$, and its closure agrees with $L_{\mathbf{E}}$. \end{theorem}

\begin{proof} Since $L_{\mathbf{E},0}$ is symmetric and nonnegative, Proposition~\ref{prop:criterio-autoadjuncion-esencial-semibounded} reduces the proof to verifying that $\ker(L_{\mathbf{E},0}^*+I)=\{0\}$. Let $\mathbf{u}$ belong to this kernel. Then \[
L_{\mathbf{E},0}\mathbf{u}=-\mathbf{u}
\] in the distributional sense. Theorem~\ref{teo:regularidad-eliptica-local-operadores-haces}, applied to the elliptic operator $L_{\mathbf{E},0}+I$, implies that $\mathbf{u}$ is smooth.

Take the cutoff functions $(\chi_j)$ from Lemma~\ref{lem:funciones-corte-completitud-capitulo-calor}. The Leibniz identities \[
\nabla^{\mathbf{E}}(\chi_j\mathbf{u})=d\chi_j\otimes \mathbf{u}+\chi_j\nabla^{\mathbf{E}}\mathbf{u}
\] and \[
\nabla^{\mathbf{E}}(\chi_j^2\mathbf{u})=2\chi_jd\chi_j\otimes \mathbf{u}+\chi_j^2\nabla^{\mathbf{E}}\mathbf{u}
\] give, after expanding and taking real parts, \[
|\nabla^{\mathbf{E}}(\chi_j\mathbf{u})|_{\mathbf{g},\mathbf{h}_{\mathbf{E}}}^2
=
\operatorname{Re}\langle\nabla^{\mathbf{E}}\mathbf{u},\nabla^{\mathbf{E}}(\chi_j^2\mathbf{u})\rangle_{\mathbf{g},\mathbf{h}_{\mathbf{E}}}
+|d\chi_j\otimes \mathbf{u}|_{\mathbf{g},\mathbf{h}_{\mathbf{E}}}^2.
\] The section $\chi_j^2\mathbf{u}$ is smooth and compactly supported. Integrating and using the definition of the distributional adjoint, \[
\|\nabla^{\mathbf{E}}(\chi_j\mathbf{u})\|_{L^2(M,T^*M\otimes \mathbf{E})}^2
=
\operatorname{Re}\langle L_{\mathbf{E},0}^*\mathbf{u},\chi_j^2\mathbf{u}\rangle_{L^2(M,\mathbf{E})}
+
\|d\chi_j\otimes \mathbf{u}\|_{L^2(M,T^*M\otimes \mathbf{E})}^2.
\] Since $L_{\mathbf{E},0}^*\mathbf{u}=-\mathbf{u}$, \[
0\leq
\|\nabla^{\mathbf{E}}(\chi_j\mathbf{u})\|_{L^2(M,T^*M\otimes \mathbf{E})}^2
=
-\|\chi_j\mathbf{u}\|_{L^2(M,\mathbf{E})}^2
+
\|d\chi_j\otimes \mathbf{u}\|_{L^2(M,T^*M\otimes \mathbf{E})}^2.
\] Therefore \[
\|\chi_j\mathbf{u}\|_{L^2(M,\mathbf{E})}^2
\leq
\|d\chi_j\|_{L^\infty(M,T^*M)}^2\|\mathbf{u}\|_{L^2(M,\mathbf{E})}^2
\longrightarrow0.
\] On the other hand, $|\chi_j\mathbf{u}-\mathbf{u}|_{\mathbf{h}_{\mathbf{E}}}^2=|1-\chi_j|^2|\mathbf{u}|_{\mathbf{h}_{\mathbf{E}}}^2\leq|\mathbf{u}|_{\mathbf{h}_{\mathbf{E}}}^2$, the dominating function belongs to $L^1(M)$, and $\chi_j\mathbf{u}\to \mathbf{u}$ pointwise. Lebesgue's dominated convergence theorem~\ref{convergencia dominada} then gives $\|\chi_j\mathbf{u}-\mathbf{u}\|_{L^2(M,\mathbf{E})}\to0$. Since at the same time $\|\chi_j\mathbf{u}\|_{L^2(M,\mathbf{E})}\to0$, it follows that $\mathbf{u}=0$. \end{proof}

\subsection{The spectral semigroup and the scalar semigroup}

For $t\geq0$, define \[
P_t^{\mathbf{E}}:=e^{-tL_{\mathbf{E}}}
\] by spectral calculus. On the trivial bundle of rank one, we write $L=\Delta_B=-\Delta$ for the positive scalar realization and $P_t=e^{-tL}$.

\begin{theorem}[Covariant heat semigroup on $L^2(M,\mathbf E)$] \label{teo:semigrupo-calor-covariante-L2} Let $(M,\mathbf{g})$ be a connected, complete Riemannian manifold without boundary, and let $\mathbf{E}\to M$ be a smooth vector bundle of finite rank with a bundle metric and a compatible connection. The family $(P_t^{\mathbf{E}})_{t\geq0}$ satisfies: \begin{enumerate}[label=(\alph*)] \item $P_0^{\mathbf{E}}=I$ and $P_{t+s}^{\mathbf{E}}=P_t^{\mathbf{E}}P_s^{\mathbf{E}}$ for $s,t\geq0$; \item each $P_t^{\mathbf{E}}$ is self-adjoint and $\|P_t^{\mathbf{E}}\|_{\mathcal L(L^2(M,\mathbf{E}))}\leq1$; \item $P_t^{\mathbf{E}}\mathbf{u}\to \mathbf{u}$ in $L^2(M,\mathbf{E})$ as $t\to0^+$; \item if $t>0$, then $P_t^{\mathbf{E}}\mathbf{u}\in\mathcal D(L_{\mathbf{E}}^m)$ for every $m\in\mathbb N$; \item for $\mathbf{u}\in L^2(M,\mathbf{E})$, the map $t\mapsto P_t^{\mathbf{E}}\mathbf{u}$ is differentiable in $(0,\infty)$ and \[
\frac{d}{dt}P_t^{\mathbf{E}}\mathbf{u}=-L_{\mathbf{E}}P_t^{\mathbf{E}}\mathbf{u};
\] if $\mathbf{u}\in\mathcal D(L_{\mathbf{E}})$, the right derivative at $t=0$ exists and equals $-L_{\mathbf{E}}\mathbf{u}$. \end{enumerate} \end{theorem}

\begin{proof} This is Proposition~\ref{prop:semigrupo-espectral-autoadjunto} applied to $A=L_{\mathbf{E}}$. For regularization, it suffices to observe that, if $t>0$, the function $\lambda\mapsto\lambda^me^{-t\lambda}$ is bounded on $[0,\infty)$. \end{proof}

Before comparing the scalar and covariant semigroups, we need positivity of the scalar resolvent.

\begin{lemma}[Positivity and uniform bound for the scalar resolvent] \label{lem:resolvente-escalar-submarkoviano} Let $(M,\mathbf{g})$ be a connected, complete Riemannian manifold without boundary. Let $\lambda>0$ and let $R_\lambda=(\lambda I+L)^{-1}$. Then: \begin{enumerate}[label=(\alph*)] \item if $f\in L^2(M)$ is real and $f\geq0$, then $R_\lambda f\geq0$ almost everywhere; \item if $f\in L^2(M)\cap L^\infty(M)$, then \[
\|R_\lambda f\|_{L^\infty(M)}
\leq
\frac{1}{\lambda}\|f\|_{L^\infty(M)}.
\] \end{enumerate} \end{lemma}

\begin{proof} Let $u=R_\lambda f$. The variational characterization of the resolvent, Corollary~\ref{cor:caracterizacion-variacional-resolvente-forma}, gives \[
\int_M\langle du,d\varphi\rangle_{\mathbf{g}}\,d\lambda_{\mathbf{g}}
+
\lambda\int_Mu\varphi\,d\lambda_{\mathbf{g}}
=
\int_Mf\varphi\,d\lambda_{\mathbf{g}}
\] for every $\varphi\in W^{1,2}(M)$.

First suppose that $f\geq0$ and take $\displaystyle u_-:=\displaystyle\max\{-u,0\}$. Theorem~\ref{teo:regla-cadena-lipschitz-sobolev}, applied to the negative part, gives $u_-\in W^{1,2}(M)$ and \[
du_-=-\mathbf 1_{\{u<0\}}du
\] almost everywhere. Choosing $\varphi=u_-$ yields \[
-\|du_-\|_{L^2(M,T^*M)}^2
-\lambda\|u_-\|_{L^2(M)}^2
=
\int_Mfu_-\,d\lambda_{\mathbf{g}}
\geq0.
\] The left-hand side is nonpositive; therefore $u_-=0$.

For the uniform bound, we may first assume that $f$ is real. The resolvent preserves real-valued functions because its defining form has real coefficients. Let $c=\displaystyle\frac{\|f\|_{L^\infty(M)}}{\lambda}$. If $c=0$, the first part applied to $u$ and $-u$ gives $u=0$. Thus suppose that $c>0$ and set $z=(u-c)_+$. The function $z$ belongs to $W^{1,2}(M)$; moreover, the set $\{u>c\}$ has finite measure because $u\in L^2(M)$, and the Cauchy--Schwarz inequality in Theorem~\ref{teo:cauchy-schwarz-hilbert} gives $z\in L^1(M)$. On $\{u>c\}$, we have $u=z+c$, whereas $dz=du$. Taking $\varphi=z$, \[
\|dz\|_{L^2(M,T^*M)}^2
+
\lambda\|z\|_{L^2(M)}^2
+
\lambda c\int_Mz\,d\lambda_{\mathbf{g}}
=
\int_Mfz\,d\lambda_{\mathbf{g}}
\leq
\lambda c\int_Mz\,d\lambda_{\mathbf{g}}.
\] Thus $z=0$ and $u\leq c$. Applying the same argument to $-u$ gives $|u|\leq c$. Now let $f$ be complex-valued, and write again $u=R_\lambda f$. The resolvent preserves real-valued functions and is positive. For each $\theta\in\mathbb R$, \[
\operatorname{Re}(e^{-i\theta}u)
=
R_\lambda\bigl(\operatorname{Re}(e^{-i\theta}f)\bigr)
\leq
R_\lambda|f|.
\] The identity $|z|=\displaystyle\sup_{\theta\in2\pi\mathbb Q}\operatorname{Re}(e^{-i\theta}z)$, valid for every $z\in\mathbb C$, and countability of $2\pi\mathbb Q$ allow us to take the supremum outside a single null set. Thus $|R_\lambda f|\leq R_\lambda|f|$ almost everywhere. Applying the real-valued bound already proved to $|f|$ yields $|u|\leq\displaystyle\frac{1}{\lambda}\|f\|_{L^\infty(M)}$. \end{proof}

\begin{corollary}[Sub-Markov property of the scalar semigroup] \label{cor:semigrupo-escalar-submarkoviano} Let $(M,\mathbf{g})$ be a connected, complete Riemannian manifold without boundary. The semigroup $(P_t)_{t\geq0}$ preserves positivity. Moreover, for $f\in L^2(M)\cap L^\infty(M)$, \[
\|P_tf\|_{L^\infty(M)}\leq\|f\|_{L^\infty(M)}.
\] \end{corollary}

\begin{proof} The exponential formula in Proposition~\ref{prop:formula-exponencial-semigrupo}, applied to the generator $-L$, gives \[
P_tf
=
\lim_{n\to\infty}
\left(I+\frac{t}{n}L\right)^{-n}f
\] in $L^2(M)$. Each factor is \[
\left(I+\frac{t}{n}L\right)^{-1}
=
\frac{n}{t}
\left(\frac{n}{t}I+L\right)^{-1}
\] and, by the preceding lemma, preserves positivity and is a contraction on $L^\infty(M)$. The iterates have the same properties. From convergence in $L^2(M)$, extract a subsequence converging almost everywhere; passage to the limit preserves positivity and the pointwise bound. \end{proof}

\subsection{The Kato--Simon inequality}

The weak Kato inequality, Corollary~\ref{desigualdad debil de kato general}, asserts that if $\mathbf{u}\in W^{1,2}(M,\mathbf{E})$, then $|\mathbf{u}|_{\mathbf{h}_{\mathbf{E}}}\in W^{1,2}(M)$ and \[
|d|\mathbf{u}|_{\mathbf{h}_{\mathbf{E}}}|_{\mathbf{g}}\leq|\nabla^{\mathbf{E}}\mathbf{u}|_{\mathbf{g},\mathbf{h}_{\mathbf{E}}}
\] almost everywhere. To obtain semigroup domination, we first prove a resolvent inequality.

\begin{lemma}[Resolvent domination] \label{lem:dominacion-resolventes-conexion} Let $(M,\mathbf{g})$ be a connected, complete Riemannian manifold without boundary, and let $\mathbf{E}\to M$ be a smooth vector bundle of finite rank with a bundle metric and a compatible connection. For $\lambda>0$ and $\mathbf{f}\in L^2(M,\mathbf{E})$, \[
\left|(\lambda I+L_{\mathbf{E}})^{-1}\mathbf{f}\right|_{\mathbf{h}_{\mathbf{E}}}
\leq
(\lambda I+L)^{-1}|\mathbf{f}|_{\mathbf{h}_{\mathbf{E}}}
\] almost everywhere. \end{lemma}

\begin{proof} Set \[
\mathbf{u}:=(\lambda I+L_{\mathbf{E}})^{-1}\mathbf{f}
\] and, for $\varepsilon>0$, \[
a_\varepsilon:=\bigl(|\mathbf{u}|_{\mathbf{h}_{\mathbf{E}}}^2+\varepsilon^2\bigr)^{\frac{1}{2}}.
\] The function $a_\varepsilon$ belongs to $W^{1,2}_{\mathrm{loc}}(M)$, and \[
da_\varepsilon(\mathbf{X})
=
\frac{\operatorname{Re}\langle\nabla_{\mathbf{X}}^{\mathbf{E}}\mathbf{u},\mathbf{u}\rangle_{\mathbf{h}_{\mathbf{E}}}}{a_\varepsilon}
\] for almost every point and every vector field $\mathbf{X}$. In particular, \[
|da_\varepsilon|_{\mathbf{g}}\leq|\nabla^{\mathbf{E}}\mathbf{u}|_{\mathbf{g},\mathbf{h}_{\mathbf{E}}}.
\]

Let $0\leq\varphi\in C_c^\infty(M)$. The section \[
\mathbf{v}_\varepsilon:=\frac{\varphi}{a_\varepsilon}\mathbf{u}
\] belongs to $W^{1,2}(M,\mathbf{E})$ and has compact support. The variational resolvent equation is the one in Corollary~\ref{cor:caracterizacion-variacional-resolvente-forma} applied to the form $\mathfrak q_{\mathbf{E}}$: \[
\mathfrak q_{\mathbf{E}}(\mathbf{u},\mathbf{v}_\varepsilon)
+
\lambda\langle \mathbf{u},\mathbf{v}_\varepsilon\rangle_{L^2(M,\mathbf{E})}
=
\langle \mathbf{f},\mathbf{v}_\varepsilon\rangle_{L^2(M,\mathbf{E})}.
\] The Leibniz rule gives \[
\nabla^{\mathbf{E}}\mathbf{v}_\varepsilon
=
\frac{d\varphi}{a_\varepsilon}\otimes \mathbf{u}
+
\frac{\varphi}{a_\varepsilon}\nabla^{\mathbf{E}}\mathbf{u}
-
\frac{\varphi}{a_\varepsilon^2}da_\varepsilon\otimes \mathbf{u}.
\] Contracting with $\nabla^{\mathbf{E}}\mathbf{u}$ and taking real parts, the first term produces $\langle da_\varepsilon,d\varphi\rangle_{\mathbf{g}}$. For the other two terms, we obtain \[
\operatorname{Re}\left\langle
\nabla^{\mathbf{E}}\mathbf{u},
\frac{\varphi}{a_\varepsilon}\nabla^{\mathbf{E}}\mathbf{u}
-
\frac{\varphi}{a_\varepsilon^2}da_\varepsilon\otimes \mathbf{u}
\right\rangle_{\mathbf{g},\mathbf{h}_{\mathbf{E}}}
=
\frac{\varphi}{a_\varepsilon}
\left(|\nabla^{\mathbf{E}}\mathbf{u}|_{\mathbf{g},\mathbf{h}_{\mathbf{E}}}^2-|da_\varepsilon|_{\mathbf{g}}^2\right)
\geq0.
\] Therefore \[
\operatorname{Re}\mathfrak q_{\mathbf{E}}(\mathbf{u},\mathbf{v}_\varepsilon)
\geq
\int_M\langle da_\varepsilon,d\varphi\rangle_{\mathbf{g}}\,d\lambda_{\mathbf{g}}.
\] Moreover, \[
\operatorname{Re}\langle \mathbf{u},\mathbf{v}_\varepsilon\rangle_{\mathbf{h}_{\mathbf{E}}}
=
\varphi\frac{|\mathbf{u}|_{\mathbf{h}_{\mathbf{E}}}^2}{a_\varepsilon}
=
\varphi\left(a_\varepsilon-\frac{\varepsilon^2}{a_\varepsilon}\right),
\] and \[
\operatorname{Re}\langle \mathbf{f},\mathbf{v}_\varepsilon\rangle_{\mathbf{h}_{\mathbf{E}}}
\leq
\varphi|\mathbf{f}|_{\mathbf{h}_{\mathbf{E}}}\frac{|\mathbf{u}|_{\mathbf{h}_{\mathbf{E}}}}{a_\varepsilon}
\leq
\varphi|\mathbf{f}|_{\mathbf{h}_{\mathbf{E}}}.
\] Substituting these estimates into the variational equation, \[
\int_M\langle da_\varepsilon,d\varphi\rangle_{\mathbf{g}}\,d\lambda_{\mathbf{g}}
+
\lambda\int_Ma_\varepsilon\varphi\,d\lambda_{\mathbf{g}}
\leq
\int_M|\mathbf{f}|_{\mathbf{h}_{\mathbf{E}}}\varphi\,d\lambda_{\mathbf{g}}
+
\lambda\varepsilon\int_M\varphi\,d\lambda_{\mathbf{g}},
\] since $\displaystyle\frac{\varepsilon^2}{a_\varepsilon}\leq\varepsilon$.

Now let $\varepsilon\to0^+$. We have $a_\varepsilon\to|\mathbf{u}|_{\mathbf{h}_{\mathbf{E}}}$ pointwise and $|a_\varepsilon-|\mathbf{u}|_{\mathbf{h}_{\mathbf{E}}}|\leq\varepsilon$, so the zeroth-order term converges. Likewise, \[
da_\varepsilon
=
\frac{\operatorname{Re}\langle\nabla^{\mathbf{E}}\mathbf{u},\mathbf{u}\rangle_{\mathbf{h}_{\mathbf{E}}}}{a_\varepsilon}
\longrightarrow d|\mathbf{u}|_{\mathbf{h}_{\mathbf{E}}}
\] almost everywhere; on the set where $\mathbf{u}=0$, the weak differential of $|\mathbf{u}|_{\mathbf{h}_{\mathbf{E}}}$ vanishes almost everywhere. Since $|da_\varepsilon|_{\mathbf{g}}\leq|\nabla^{\mathbf{E}}\mathbf{u}|_{\mathbf{g},\mathbf{h}_{\mathbf{E}}}$ and $d\varphi$ is bounded and compactly supported, Lebesgue's dominated convergence theorem~\ref{convergencia dominada} gives \[
\int_M\langle d|\mathbf{u}|_{\mathbf{h}_{\mathbf{E}}},d\varphi\rangle_{\mathbf{g}}\,d\lambda_{\mathbf{g}}
+
\lambda\int_M|\mathbf{u}|_{\mathbf{h}_{\mathbf{E}}}\varphi\,d\lambda_{\mathbf{g}}
\leq
\int_M|\mathbf{f}|_{\mathbf{h}_{\mathbf{E}}}\varphi\,d\lambda_{\mathbf{g}}.
\]

Let \[
w:=(\lambda I+L)^{-1}|\mathbf{f}|_{\mathbf{h}_{\mathbf{E}}}.
\] Subtracting the weak equation for $w$, we wish to use the test function \[
z:=(|\mathbf{u}|_{\mathbf{h}_{\mathbf{E}}}-w)_+.
\] Theorem~\ref{teo:regla-cadena-lipschitz-sobolev}, applied to $|\mathbf{u}|_{\mathbf{h}_{\mathbf{E}}}-w$ and the function $t\mapsto t_+$, shows that $z\in W^{1,2}(M)$ and that \[
dz=\mathbf 1_{\{|\mathbf{u}|_{\mathbf{h}_{\mathbf{E}}}>w\}}(d|\mathbf{u}|_{\mathbf{h}_{\mathbf{E}}}-dw)
\] almost everywhere. By Proposition~\ref{prop:W12-cero-igual-W12-completa}, $W^{1,2}(M)=W_0^{1,2}(M)$. To use $z$ as a test function, we combine smooth approximation with nonnegative truncations. Choose $\phi_j\in C_c^\infty(M)$ with $\phi_j\to z$ in $W^{1,2}(M)$ and, after passing to a subsequence, assume that $\phi_j\to z$ almost everywhere. For $\delta>0$, let $\vartheta_\delta\in C^\infty(\mathbb R)$ be a convex function such that \[
\vartheta_\delta(t)=0\quad\text{if }t\leq0,
\qquad
0\leq\vartheta_\delta'(t)\leq1,
\qquad
|\vartheta_\delta(t)-t_+|\leq\delta.
\] These functions are obtained by integrating a smooth nondecreasing approximation to the indicator function of $(0,+\infty)$. Choose $\delta_j\to0^+$ so that \[
\delta_j\,\lambda_{\mathbf{g}}(\operatorname{supp}\phi_j)^{\frac{1}{2}}\longrightarrow0
\] and define $z_j:=\vartheta_{\delta_j}(\phi_j)$. Then $z_j\in C_c^\infty(M)$, $z_j\geq0$, and \[
\|z_j-z\|_{L^2(M)}
\leq
\|\phi_j-z\|_{L^2(M)}
+
\delta_j\lambda_{\mathbf{g}}(\operatorname{supp}\phi_j)^{\frac{1}{2}}
\longrightarrow0.
\] The chain rule gives $dz_j=\vartheta_{\delta_j}'(\phi_j)d\phi_j$. Therefore \[
|dz_j-dz|_{\mathbf{g}}
\leq
|d\phi_j-dz|_{\mathbf{g}}
+
|\vartheta_{\delta_j}'(\phi_j)-\mathbf 1_{\{z>0\}}|\,|dz|_{\mathbf{g}}.
\] The first term converges to zero in $L^2(M,T^*M)$. In the second, the factor multiplying $|dz|_{\mathbf{g}}$ converges to zero almost everywhere: this is immediate on $\{z>0\}$, whereas $dz=0$ almost everywhere on the level set $\{z=0\}$. Since this factor is bounded by $1$, Lebesgue's dominated convergence theorem~\ref{convergencia dominada} shows that $dz_j\to dz$ in $L^2(M,T^*M)$. Thus $z_j\to z$ in $W^{1,2}(M)$.

The preceding weak inequality is continuous with respect to this norm, since $d|\mathbf{u}|_{\mathbf{h}_{\mathbf{E}}}-dw\in L^2(M,T^*M)$ and $\lambda(|\mathbf{u}|_{\mathbf{h}_{\mathbf{E}}}-w)\in L^2(M)$. We may therefore apply it to $z_j$ and pass to the limit. On the set $\{z>0\}$, we have $|\mathbf{u}|_{\mathbf{h}_{\mathbf{E}}}-w=z$, and obtain \[
\int_M|dz|_{\mathbf{g}}^2\,d\lambda_{\mathbf{g}}
+
\lambda\int_Mz^2\,d\lambda_{\mathbf{g}}
\leq0.
\] Both terms are nonnegative, so $z=0$ almost everywhere. This is equivalent to $|\mathbf{u}|_{\mathbf{h}_{\mathbf{E}}}\leq w$. \end{proof}

\begin{theorem}[Kato--Simon inequality] \label{teo:kato-simon-laplaciano-conexion} Let $(M,\mathbf{g})$ be a connected, complete Riemannian manifold without boundary, and let $\mathbf{E}\to M$ be a smooth vector bundle of finite rank with a bundle metric and a compatible connection. \index{Kato--Simon inequality} For every $t\geq0$ and every $\mathbf{u}\in L^2(M,\mathbf{E})$, \[
|P_t^{\mathbf{E}}\mathbf{u}|_{\mathbf{h}_{\mathbf{E}}}\leq P_t(|\mathbf{u}|_{\mathbf{h}_{\mathbf{E}}})
\] almost everywhere. \end{theorem}

\begin{proof} For $\tau>0$, the preceding lemma with $\lambda=\displaystyle\frac{1}{\tau}$ is equivalent to \[
\left|(I+\tau L_{\mathbf{E}})^{-1}v\right|_{\mathbf{h}_{\mathbf{E}}}
\leq
(I+\tau L)^{-1}|v|_{\mathbf{h}_{\mathbf{E}}}.
\] Fix $\tau>0$ and write $R_{\mathbf E}=(I+\tau L_{\mathbf E})^{-1}$ and $R=(I+\tau L)^{-1}$. The preceding inequality is the case $k=1$ of $|R_{\mathbf E}^{k}\mathbf u|\leq R^{k}|\mathbf u|$. If it holds for $k\in\mathbb N$, apply it first to $R_{\mathbf E}^{k}\mathbf u$ and then use the fact that $R$ preserves positivity: \[
 |R_{\mathbf E}^{k+1}\mathbf u|
 \leq R\bigl(|R_{\mathbf E}^{k}\mathbf u|\bigr)
 \leq R\bigl(R^{k}|\mathbf u|\bigr)
 =R^{k+1}|\mathbf u|.
\] This proves the inequality for each $k\in\mathbb N$. Given $t>0$, take $\tau=\displaystyle\frac{t}{N}$ and $k=N$, with $N\in\mathbb N$, to obtain \[
 \left|\left(I+\frac{t}{N}L_{\mathbf E}\right)^{-N}\mathbf u\right|_{\mathbf h_{\mathbf E}}
 \leq\left(I+\frac{t}{N}L\right)^{-N}|\mathbf u|_{\mathbf h_{\mathbf E}}.
\] Proposition~\ref{prop:formula-exponencial-semigrupo} shows that the left-hand side converges in $L^2(M,\mathbf{E})$ to $P_t^{\mathbf{E}}\mathbf{u}$ and the right-hand side in $L^2(M)$ to $P_t|\mathbf{u}|_{\mathbf{h}_{\mathbf{E}}}$. Choose a common subsequence converging almost everywhere and pass to the limit in the pointwise inequality. The resulting inequality is the assertion for $t>0$; at $t=0$, both sides agree. \end{proof}

\section{The covariant heat kernel}

The kernel $\mathbf{K}_t^{\mathbf{E}}$ is a section of $\mathbf{E}\boxtimes \mathbf{E}^*\to M_x\times M_y$, equipped with the product connection. In every differential calculation, we distinguish the variables by writing \[
\nabla_x^r\nabla_y^s\mathbf{K}_t^{\mathbf{E}}.
\] After fixing charts or frames, we write its components as \[
\bigl(\nabla_x^r\nabla_y^s\mathbf{K}_t^{\mathbf{E}}\bigr)_{
i_1\cdots i_r\,;\,j_1\cdots j_s},
\] where the semicolon separates the indices of $x$ from those of $y$. Each group follows the recursive convention in Chapter~\ref{cap:sobolev-haces}; in particular, these derivatives will not be mixed with undeclared indexed compositions. We denote by \[
\Delta_{B,\mathbf{E}}^{(x)}
\quad\text{and}\quad
\Delta_{B,\mathbf{E}^*}^{(y)}
\] the positive Bochner Laplacians acting, respectively, on the factor $\mathbf{E}_x$ and the dual factor $\mathbf{E}_y^*$ of the product connection.

\subsection{Regularization and continuity of evaluations}

The essential step in constructing the kernel is that $P_t^{\mathbf{E}}$ transforms sections in $L^2(M,\mathbf{E})$ into smooth sections when $t>0$. This property follows directly from spectral calculus and the elliptic regularity developed in the preceding chapter.

\begin{lemma}[Semigroup regularization] \label{lem:regularizacion-semigrupo-calor-covariante} Let $(M,\mathbf{g})$ be a connected, complete Riemannian manifold without boundary, and let $\mathbf{E}\to M$ be a smooth vector bundle of finite rank with a bundle metric and a compatible connection. Let $t>0$, $k\in\mathbb N_0$, and let $K\subseteq M$ be compact. The operator $P_t^{\mathbf{E}}$ induces a continuous linear map \[
P_t^{\mathbf{E}}\colon L^2(M,\mathbf{E})\longrightarrow C^k(K,\mathbf{E}).
\] In particular, $P_t^{\mathbf{E}}\mathbf{u}$ has a smooth representative for every $\mathbf{u}\in L^2(M,\mathbf{E})$. \end{lemma}

\begin{proof}
For $\mathbf w\in\mathcal D(L_{\mathbf E})$, the definition by closure gives sections $\mathbf w_j\in\Gamma_c(\mathbf E)$ such that $\mathbf w_j\to\mathbf w$ and $L_{\mathbf E,0}\mathbf w_j\to L_{\mathbf E}\mathbf w$ in $L^2$. Distributional continuity of the differential operator implies $L_{\mathbf E,0}\mathbf w=L_{\mathbf E}\mathbf w$ distributionally. Spectral calculus shows that, if $\mathbf v\in\mathcal D(L_{\mathbf E}^m)$, then $L_{\mathbf E}^j\mathbf v\in\mathcal D(L_{\mathbf E})$ for $j<m$. Successive applications of the preceding identity give
\[
L_{\mathbf E,0}^j\mathbf v=L_{\mathbf E}^j\mathbf v\in L^2(M,\mathbf E),
\qquad 0\leq j\leq m.
\]

These identities imply $\mathbf v\in H^{2m}_{\mathrm{loc}}(M,\mathbf E)$. For $m=1$ this is Theorem~\ref{teo:regularidad-eliptica-local-operadores-haces}. If the assertion holds through $m$, apply it to $L_{\mathbf E,0}\mathbf v$ when powers through $m+1$ exist; thus $L_{\mathbf E,0}\mathbf v\in H^{2m}_{\mathrm{loc}}$. The same theorem, with right-hand side in that space, gives $\mathbf v\in H^{2m+2}_{\mathrm{loc}}$. In particular, $P_t^{\mathbf E}\mathbf u$, which belongs to every domain by spectral calculus, has a smooth representative by local Sobolev embeddings.

Justify the regularization bound without assuming uniform geometry. Fix $m\geq1$, a trivializing chart, and a smooth cutoff $\chi$ compactly supported in it. The map assigning to $\mathbf v\in\mathcal D(L_{\mathbf E}^m)$ the components of $\chi\mathbf v$, extended by zero, takes values in $H^{2m}(\mathbb R^n,\mathbb K^r)$ by the local regularity already proved. Its graph is closed: if $\mathbf v_j\to\mathbf v$ in graph norm and the components of $\chi\mathbf v_j$ converge to $w$ in $H^{2m}$, both convergences imply $L^2$ convergence and force $w$ to be the components of $\chi\mathbf v$. The domain with its graph norm is Banach because $L_{\mathbf E}^m$ is closed. Theorem~\ref{teo:grafica-cerrada} therefore gives
\[
\|\chi\mathbf v\|_{H^{2m}}
\leq C_\chi\bigl(\|\mathbf v\|_{L^2(M,\mathbf E)}
                       +\|L_{\mathbf E}^m\mathbf v\|_{L^2(M,\mathbf E)}\bigr),
\]
where the left-hand norm is taken in coordinate components. Choose finitely many such cutoffs equal to one on open sets covering $K$. Sobolev embedding, with $2m>k+\displaystyle\frac{n}{2}$, and norm equivalence on these compact sets give
\[
\|\mathbf v\|_{C^k(K,\mathbf E)}
\leq C_{K,k,m}\bigl(\|\mathbf v\|_{L^2(M,\mathbf E)}
                      +\|L_{\mathbf E}^m\mathbf v\|_{L^2(M,\mathbf E)}\bigr).
\]
By spectral calculus,
\[
\|L_{\mathbf E}^mP_t^{\mathbf E}\|_{\mathcal L(L^2(M,\mathbf E))}
\leq\sup_{\lambda\geq0}\lambda^m e^{-t\lambda}
=\left(\frac{m}{et}\right)^m.
\]
Together with contractivity of $P_t^{\mathbf E}$, this yields
\begin{equation}
\label{eq:regularizacion-calor-cota-tiempo-compacto}
\|P_t^{\mathbf E}\mathbf u\|_{C^k(K,\mathbf E)}
\leq C_{K,k,m}\left(1+\left(\frac{m}{et}\right)^m\right)
                         \|\mathbf u\|_{L^2(M,\mathbf E)}.
\end{equation}
The bound is uniform for $t\geq t_0>0$. Its constant may depend on $K$; completeness alone does not provide uniformity with respect to its position.
\end{proof}

For $t>0$ and $x\in M$, denote the regularized evaluation operator by \[
R_{t,x}\colon L^2(M,\mathbf{E})\longrightarrow \mathbf{E}_x,
\qquad
R_{t,x}\mathbf{u}:=P_t^{\mathbf{E}}\mathbf{u}(x),
\]. The preceding lemma shows that $R_{t,x}$ is continuous.

\begin{theorem}[Covariant heat kernel] \label{teo:nucleo-calor-covariante-suave} Let $(M,\mathbf{g})$ be a connected, complete Riemannian manifold without boundary, and let $\mathbf{E}\to M$ be a smooth vector bundle of finite rank with a bundle metric and a compatible connection. \index{heat kernel@heat kernel!covariant} There exists a unique smooth map \[
\mathbf{K}^{\mathbf{E}}\colon (0,\infty)\times M\times M\longrightarrow \mathbf{E}\boxtimes \mathbf{E}^*,
\qquad
(t,x,y)\longmapsto \mathbf{K}_t^{\mathbf{E}}(x,y)\in\operatorname{Hom}(\mathbf{E}_y,\mathbf{E}_x),
\] such that, for every $t>0$, $\mathbf{u}\in L^2(M,\mathbf{E})$, and $x\in M$, \[
P_t^{\mathbf{E}}\mathbf{u}(x)
=
\int_M\mathbf{K}_t^{\mathbf{E}}(x,y)\mathbf{u}(y)\,d\lambda_{\mathbf{g}}(y).
\] The integral converges absolutely in the finite-dimensional fiber $\mathbf{E}_x$. Moreover: \begin{enumerate}[label=(\alph*)] \item $\mathbf{K}_t^{\mathbf{E}}(y,x)=\mathbf{K}_t^{\mathbf{E}}(x,y)^*$; \item for $s,t>0$, \[
\mathbf{K}_{t+s}^{\mathbf{E}}(x,y)
=
\int_M\mathbf{K}_t^{\mathbf{E}}(x,z)\mathbf{K}_s^{\mathbf{E}}(z,y)\,d\lambda_{\mathbf{g}}(z);
\] \item the kernel satisfies, smoothly in both variables, \[
\frac{\partial}{\partial t}\mathbf{K}_t^{\mathbf{E}}
=-\Delta_{B,\mathbf{E}}^{(x)}\mathbf{K}_t^{\mathbf{E}}
=-\Delta_{B,\mathbf{E}^*}^{(y)}\mathbf{K}_t^{\mathbf{E}}.
\] In particular, for $y\in M$ and $\boldsymbol{\eta}\in \mathbf{E}_y$, the first equality follows by evaluating the kernel on $\boldsymbol{\eta}$; \item if $p_t(x,y)$ is the kernel obtained by the same construction for the trivial bundle of rank one, then \[
\|\mathbf{K}_t^{\mathbf{E}}(x,y)\|_{\mathcal L((\mathbf{E}_y,\mathbf{h}_{\mathbf{E}}(y)),(\mathbf{E}_x,\mathbf{h}_{\mathbf{E}}(x)))}
\leq
p_t(x,y).
\] \end{enumerate} \end{theorem}

\begin{proof} Fix $t>0$ and set $s=\displaystyle\frac{t}{2}$. Since $R_{s,x}$ is continuous, it has an adjoint \[
R_{s,x}^*\colon \mathbf{E}_x\longrightarrow L^2(M,\mathbf{E}).
\] Define \[
\mathbf{K}_t^{\mathbf{E}}(x,y):=R_{s,x}R_{s,y}^*.
\] This formula already shows that $\mathbf{K}_t^{\mathbf{E}}(x,y)$ is a linear map from $\mathbf{E}_y$ to $\mathbf{E}_x$.

\emph{Smoothness.} Let $U,V\subseteq M$ be open subsets over which smooth orthonormal frames $(\mathbf{e}_a)_{a=1}^r$ and $(\mathbf{f}_b)_{b=1}^r$ of $\mathbf{E}$ exist. For $s>0$ and $y\in V$, define \[
\boldsymbol{\Psi}_b(s,y):=R_{s,y}^*\mathbf{f}_b(y)\in L^2(M,\mathbf{E}).
\] Fix $\mathbf{u}\in L^2(M,\mathbf{E})$. By the definition of the adjoint, \begin{equation}
\label{eq:componentes-debiles-Psi-nucleo-calor}
\langle \mathbf{u},\boldsymbol{\Psi}_b(s,y)\rangle_{L^2(M,\mathbf{E})}
=
\langle P_s^{\mathbf{E}}\mathbf{u}(y),\mathbf{f}_b(y)\rangle_{\mathbf{h}_{\mathbf{E}}}.
\end{equation} Let us check the joint regularity that we will use. Let $(s_0,y_0)\in(0,\infty)\times V$ and choose a compact interval $J\Subset(0,\infty)$ containing $s_0$ in its interior. Fix $\delta>0$ such that $2\delta<\inf_{s\in J}s$. For $b\in\mathbb N_0$, set
\[
C_{0,\delta}:=1,\qquad
C_{b,\delta}:=\left(\frac{b}{e\delta}\right)^b\quad(b\geq1).
\]
Maximizing $\lambda^b e^{-\delta\lambda}$ over $\lambda\geq0$ gives $\lambda^b e^{-\delta\lambda}\leq C_{b,\delta}$. If $s\in J$ and $0<|h|<\delta$, the entire segment between $s-\delta$ and $s+h-\delta$ lies in $[\delta,\infty)$. The fundamental theorem of calculus and its iterated application to the exponential give, for $\lambda\geq0$,
\[
\begin{split}
|e^{-(s+h-\delta)\lambda}-e^{-(s-\delta)\lambda}|
&\leq |h|\lambda e^{-\delta\lambda},\\
\left|\frac{e^{-(s+h-\delta)\lambda}-e^{-(s-\delta)\lambda}}{h}
       +\lambda e^{-(s-\delta)\lambda}\right|
&\leq\frac{|h|}{2}\lambda^2e^{-\delta\lambda}.
\end{split}
\]
Write $F_a(s):=(-L_{\mathbf E})^aP_{s-\delta}^{\mathbf E}\mathbf u$ for $a\in\mathbb N_0$. Spectral calculus gives
\[
\begin{split}
\|F_a(s+h)-F_a(s)\|_{L^2}
&\leq |h|C_{a+1,\delta}\|\mathbf u\|_{L^2},\\
\left\|\frac{F_a(s+h)-F_a(s)}{h}-F_{a+1}(s)\right\|_{L^2}
&\leq\frac{|h|}{2}C_{a+2,\delta}\|\mathbf u\|_{L^2}.
\end{split}
\]
In particular, both estimates are uniform in $s\in J$, and their right-hand sides tend to zero with $h$.

By Lemma~\ref{lem:regularizacion-semigrupo-calor-covariante}, for each compact set $K\Subset V$ and each $q\in\mathbb N_0$ there is $B_{K,q,\delta}>0$ such that
\[
\|P_\delta^{\mathbf E}\mathbf v\|_{C^q(K,\mathbf E)}
\leq B_{K,q,\delta}\|\mathbf v\|_{L^2(M,\mathbf E)}.
\]
Since $P_\delta^{\mathbf E}F_a(s)=(-L_{\mathbf E})^aP_s^{\mathbf E}\mathbf u$, we obtain
\[
\begin{split}
\|(-L_{\mathbf E})^a(P_{s+h}^{\mathbf E}-P_s^{\mathbf E})\mathbf u\|_{C^q(K,\mathbf E)}
&\leq B_{K,q,\delta}|h|C_{a+1,\delta}\|\mathbf u\|_{L^2},\\
\left\|\frac{(-L_{\mathbf E})^a(P_{s+h}^{\mathbf E}-P_s^{\mathbf E})\mathbf u}{h}
       -(-L_{\mathbf E})^{a+1}P_s^{\mathbf E}\mathbf u\right\|_{C^q(K,\mathbf E)}
&\leq B_{K,q,\delta}\frac{|h|}{2}C_{a+2,\delta}\|\mathbf u\|_{L^2}.
\end{split}
\]
These bounds prove by induction that $s\mapsto P_s^{\mathbf E}\mathbf u$ is smooth with values in each $C^q(K,\mathbf E)$ and that its derivative of order $a$ is $(-L_{\mathbf E})^aP_s^{\mathbf E}\mathbf u$. In a trivialization over $V$, the Leibniz rule expresses
\[
\partial_s^a\partial_y^\beta
\langle P_s^{\mathbf E}\mathbf u(y),\mathbf f_b(y)\rangle_{\mathbf h_{\mathbf E}}
\]
as a finite sum of pairings between coordinate derivatives of $(-L_{\mathbf E})^aP_s^{\mathbf E}\mathbf u$ and derivatives of the frame $\mathbf f_b$. Continuity in $C^q(K,\mathbf E)$, with $q\geq|\beta|$, and continuity of spatial derivatives imply that each summand is continuous in $(s,y)$. Consequently, the right-hand side of \eqref{eq:componentes-debiles-Psi-nucleo-calor} is jointly smooth in $(s,y)$.

We have proved that $(s,y)\mapsto\boldsymbol\Psi_b(s,y)$ is weakly smooth with values in $L^2(M,\mathbf E)$. Theorem~\ref{teo:diferenciabilidad-debil-fuerte-hilbert}, applied in charts on $(0,\infty)\times V$, shows that this map is smooth in the norm of that Hilbert space.

If $s=\displaystyle\frac{t}{2}$, then \[
\left\langle \mathbf{K}_t^{\mathbf{E}}(x,y)\mathbf{f}_b(y),\mathbf{e}_a(x)\right\rangle_{\mathbf{h}_{\mathbf{E}}}
=
\left\langle\boldsymbol{\Psi}_b(s,y),\boldsymbol{\Psi}_a(s,x)\right\rangle_{L^2(M,\mathbf{E})}.
\] The inner product is a continuous sesquilinear map; hence the last expression is smooth in $(t,x,y)$. Since this holds for every component in local frames, $\mathbf{K}^{\mathbf{E}}$ is a smooth section of $\mathbf{E}\boxtimes \mathbf{E}^*$.

\emph{Integral representation.} Fix $x\in M$ and $\xi\in \mathbf{E}_x$. For $\mathbf{u}\in L^2(M,\mathbf{E})$, the section \[
y\longmapsto \mathbf{K}_t^{\mathbf{E}}(x,y)^*\xi
=
R_{s,y}R_{s,x}^*\xi
=
P_s^{\mathbf{E}}(R_{s,x}^*\xi)(y)
\] belongs to $L^2(M,\mathbf{E})$. The Cauchy--Schwarz inequality in Theorem~\ref{teo:cauchy-schwarz-hilbert} shows that its scalar product with $\mathbf{u}$ is integrable. Using the definition of the adjoint and self-adjointness of $P_s^{\mathbf{E}}$, \begin{align*}
\int_M\left\langle \mathbf{K}_t^{\mathbf{E}}(x,y)\mathbf{u}(y),\xi\right\rangle_{\mathbf{h}_{\mathbf{E}}}\,d\lambda_{\mathbf{g}}(y)
&=
\left\langle \mathbf{u},P_s^{\mathbf{E}}R_{s,x}^*\xi\right\rangle_{L^2(M,\mathbf{E})}\\
&=
\left\langle P_s^{\mathbf{E}}\mathbf{u},R_{s,x}^*\xi\right\rangle_{L^2(M,\mathbf{E})}\\
&=
\left\langle R_{s,x}P_s^{\mathbf{E}}\mathbf{u},\xi\right\rangle_{\mathbf{h}_{\mathbf{E}}}\\
&=
\left\langle P_t^{\mathbf{E}}\mathbf{u}(x),\xi\right\rangle_{\mathbf{h}_{\mathbf{E}}}.
\end{align*} Convergence is also absolute in the fiber norm. If $(\boldsymbol\xi_a)_{a=1}^{r_{\mathbf E}}$ is an orthonormal basis of $\mathbf E_x$, then \[
 \int_M|\mathbf K_t^{\mathbf E}(x,y)\mathbf u(y)|\,d\lambda_{\mathbf g}(y)
 \leq\sum_{a=1}^{r_{\mathbf E}}
 \|\mathbf K_t^{\mathbf E}(x,\cdot)^*\boldsymbol\xi_a\|_{L^2(M,\mathbf E)}
 \|\mathbf u\|_{L^2(M,\mathbf E)}<\infty.
\] Since $\xi\in\mathbf E_x$ was arbitrary, the integral formula follows.

If another smooth kernel satisfies the same formula, their difference, paired with local frames, integrates to zero against every compactly supported smooth section. The fundamental lemma of distribution theory implies that each component vanishes; hence the kernel is unique.

\emph{Symmetry and the semigroup property.} The definition gives \[
\mathbf{K}_t^{\mathbf{E}}(y,x)^*
=
(R_{s,y}R_{s,x}^*)^*
=
R_{s,x}R_{s,y}^*
=
\mathbf{K}_t^{\mathbf{E}}(x,y).
\] For the semigroup property, fix $s,t>0$, $y\in M$, and $\boldsymbol\eta\in\mathbf E_y$. If $0<b<a$, the semigroup identity and evaluation of the resulting smooth sections give
\[
R_{a,y}=R_{b,y}P_{a-b}^{\mathbf E},\qquad
R_{a,y}^*=P_{a-b}^{\mathbf E}R_{b,y}^*.
\]
Choose $0<\varepsilon<\displaystyle\frac{s}{2}$. By the definition of the kernel and these identities,
\[
\begin{split}
\mathbf K_s^{\mathbf E}(\,\cdot\,,y)\boldsymbol\eta
&=P_{\frac{s}{2}}^{\mathbf E}R_{\frac{s}{2},y}^*\boldsymbol\eta
=P_{s-\varepsilon}^{\mathbf E}R_{\varepsilon,y}^*\boldsymbol\eta,\\
\mathbf K_{t+s}^{\mathbf E}(\,\cdot\,,y)\boldsymbol\eta
&=P_{t+s-\varepsilon}^{\mathbf E}R_{\varepsilon,y}^*\boldsymbol\eta.
\end{split}
\]
Both columns belong to $L^2(M,\mathbf E)$, and the semigroup property in that space implies
$P_t^{\mathbf E}(\mathbf K_s^{\mathbf E}(\,\cdot\,,y)\boldsymbol\eta)
=\mathbf K_{t+s}^{\mathbf E}(\,\cdot\,,y)\boldsymbol\eta$.
This equality holds pointwise because both sides are smooth. Applying the integral formula, \[
\int_M\mathbf{K}_t^{\mathbf{E}}(x,z)\mathbf{K}_s^{\mathbf{E}}(z,y)\boldsymbol{\eta}\,d\lambda_{\mathbf{g}}(z)
=
P_t^{\mathbf{E}}\bigl(\mathbf{K}_s^{\mathbf{E}}(\,\cdot\,,y)\boldsymbol{\eta}\bigr)(x)
=
\mathbf{K}_{t+s}^{\mathbf{E}}(x,y)\boldsymbol{\eta}.
\]

\emph{The heat equation.} Fix $t>0$, $y\in M$, and $\boldsymbol{\eta}\in \mathbf{E}_y$. Choose $t_0\in(0,t)$ and set \[
\mathbf{w}:=\mathbf{K}_{t_0}^{\mathbf{E}}(\,\cdot\,,y)\boldsymbol{\eta}\in L^2(M,\mathbf{E}).
\] The semigroup property gives \[
\mathbf{K}_\tau^{\mathbf{E}}(\,\cdot\,,y)\boldsymbol{\eta}=P_{\tau-t_0}^{\mathbf{E}}\mathbf{w},
\qquad \tau>t_0.
\] Apply the time estimates already proved to $\tau\mapsto P_{\tau-t_0}^{\mathbf E}\mathbf w$. Choose a compact interval $J_0\Subset(t_0,\infty)$ containing $t$ in its interior and $\delta>0$ such that $2\delta<\inf_{\tau\in J_0}(\tau-t_0)$. For each $\tau\in J_0$, each compact set $K\subseteq M$, each $q\in\mathbb N_0$, and $0<|h|<\delta$, the second estimate with $a=0$ gives
\[
\begin{split}
&\left\|\frac{\mathbf K_{\tau+h}^{\mathbf E}(\,\cdot\,,y)\boldsymbol\eta
                 -\mathbf K_\tau^{\mathbf E}(\,\cdot\,,y)\boldsymbol\eta}{h}
       +L_{\mathbf E}P_{\tau-t_0}^{\mathbf E}\mathbf w\right\|_{C^q(K,\mathbf E)}\\
&\qquad\leq B_{K,q,\delta}\frac{|h|}{2}C_{2,\delta}\|\mathbf w\|_{L^2}
\longrightarrow0.
\end{split}
\]
If $K$ is not contained in one chart, take a finite cover by charts and the maximum of their regularization constants. The right-hand side tends to zero because all its constants and $\mathbf w$ are fixed. At positive times, $P_{\tau-t_0}^{\mathbf E}\mathbf w$ belongs to the domain of every power of $L_{\mathbf E}$; on its smooth representatives, the action of $L_{\mathbf E}$ agrees with the differential expression of the connection Laplacian. Thus,
\[
\frac{\partial}{\partial t}\mathbf K_t^{\mathbf E}(\,\cdot\,,y)\boldsymbol\eta
=-\Delta_{B,\mathbf E}^{(x)}\mathbf K_t^{\mathbf E}(\,\cdot\,,y)\boldsymbol\eta.
\]
Since $q$ and $K$ were arbitrary, the equality holds in $\Gamma(\mathbf E)$. To obtain the equation in the second variable, fix $x$ and consider $A(y)=\mathbf K_t^{\mathbf E}(y,x)\in\operatorname{Hom}(\mathbf E_x,\mathbf E_y)$. Compatibility of the connection with the metric, applied to the inner product of $A(y)v$ with a local section, gives
\[
\nabla_y(A^*)=(\nabla_y A)^*.
\]
Here the induced connections on the respective homomorphism bundles are used. Repeating the identity for the second covariant derivative and taking the trace with the Riemannian metric yields
\[
\Delta_{B,\mathbf E^*}^{(y)}A^*
=(\Delta_{B,\mathbf E}^{(y)}A)^*.
\]
Symmetry $\mathbf K_t^{\mathbf E}(x,y)=\mathbf K_t^{\mathbf E}(y,x)^*$ and the equation already proved in the first variable then imply
\[
\partial_t\mathbf K_t^{\mathbf E}
=-\Delta_{B,\mathbf E^*}^{(y)}\mathbf K_t^{\mathbf E}.
\]

\emph{Pointwise domination of the kernel.} Fix $y\in M$ and $\boldsymbol{\eta}\in \mathbf{E}_y$. Let $r_{\mathbf{E}}=\operatorname{rank}\mathbf{E}$. Choose an open neighborhood $V$ of $y$ with a smooth orthonormal frame $(\mathbf{e}_a)_{a=1}^{r_{\mathbf E}}$ and write $\boldsymbol{\eta}=\displaystyle\sum_{a=1}^{r_{\mathbf{E}}}\eta_a\mathbf e_a(y)$. Define the local section of constant norm, with $\eta_a\in\mathbb K$ for $a\in\{1,\ldots,r_{\mathbf E}\}$, \[
\widetilde{\boldsymbol{\eta}}(z):=\displaystyle\sum_{a=1}^{r_{\mathbf{E}}}\eta_a\mathbf e_a(z).
\] We construct the required approximate identity explicitly. Take a chart $\kappa\colon V\longrightarrow\Omega\subseteq\mathbb R^n$ with $\kappa(y)=0$ and choose $\rho>0$ so that $\overline{B}_{\mathrm{euc}}(0,\rho)\subseteq\Omega$. Let $a\in C^\infty(\Omega)$ be the positive coordinate weight of the Riemannian measure, determined by \[
\int_V \phi\,d\lambda_{\mathbf{g}}
=
\int_\Omega(\phi\circ\kappa^{-1})(z)a(z)\,d\lambda_n(z)
\] for every integrable function $\phi$ supported in $V$. Denote by $\omega$ the mollifier in Theorem~\ref{teo:molificadores-euclidianos}; for $j\in\mathbb N$, define \[
\varepsilon_j:=\frac{\rho}{j+1},
\qquad
c_j:=
\left(
\int_{\Omega}
\varepsilon_j^{-n}\omega\left(\frac{z}{\varepsilon_j}\right)a(z)\,dz
\right)^{-1},
\] where the integrand vanishes outside $\overline{B}_{\mathrm{euc}}(0,\varepsilon_j)\subseteq\Omega$. The number $c_j$ is well defined and positive. Set \[
\rho_j(q)
:=
\begin{cases}
c_j\varepsilon_j^{-n}\omega\left(\displaystyle\frac{\kappa(q)}{\varepsilon_j}\right),&q\in V,\\
0,&q\notin V.
\end{cases}
\] Since the support of the first expression is compact and contained in $V$, its extension by zero is smooth. Moreover, $\rho_j\geq0$, $\displaystyle\int_M\rho_j\,d\lambda_{\mathbf{g}}=1$, and $\operatorname{supp}(\rho_j)\subseteq
\kappa^{-1}(\overline{B}_{\mathrm{euc}}(0,\varepsilon_j))$. If $W$ is a neighborhood of $y$, choose $0<r_W\leq\rho$ such that $\kappa^{-1}(B_{\mathrm{euc}}(0,r_W))\subseteq W$; then $\operatorname{supp}(\rho_j)\subseteq W$ for every $j$ satisfying $j+1>\displaystyle\frac{\rho}{r_W}$. Set $\mathbf{u}_j=\rho_j\widetilde{\boldsymbol{\eta}}$. Then \[
|\mathbf{u}_j|_{\mathbf{h}_{\mathbf{E}}}=|\boldsymbol{\eta}|_{\mathbf{h}_{\mathbf{E}}}\rho_j.
\] The Kato--Simon inequality initially gives, for almost every $x\in M$, \[
|P_t^{\mathbf{E}}\mathbf{u}_j(x)|_{\mathbf{h}_{\mathbf{E}}}
\leq
|\boldsymbol{\eta}|_{\mathbf{h}_{\mathbf{E}}}P_t\rho_j(x).
\] For $t>0$, both functions in this inequality are continuous, since the scalar and covariant semigroups are regularizing. If the inequality failed at a point, it would fail on a nonempty open subset and hence on a set of positive measure. It therefore holds for every $x\in M$. By the integral formula, \[
P_t^{\mathbf{E}}\mathbf{u}_j(x)
=
\int_M\mathbf{K}_t^{\mathbf{E}}(x,z)\widetilde{\boldsymbol{\eta}}(z)\rho_j(z)\,d\lambda_{\mathbf{g}}(z).
\] The function $z\mapsto \mathbf{K}_t^{\mathbf{E}}(x,z)\widetilde{\boldsymbol{\eta}}(z)$ is continuous. Given $\delta>0$, there exists a neighborhood $W$ of $y$ on which its difference from $\mathbf{K}_t^{\mathbf{E}}(x,y)\boldsymbol{\eta}$ has norm less than $\delta$ throughout $W$. The explicit choice above gives an index beyond which $\operatorname{supp}(\rho_j)\subseteq W$, and then \[
\left|
P_t^{\mathbf{E}}\mathbf{u}_j(x)-\mathbf{K}_t^{\mathbf{E}}(x,y)\boldsymbol{\eta}
\right|_{\mathbf{h}_{\mathbf{E}}}
\leq
\int_M\delta\rho_j\,d\lambda_{\mathbf{g}}
=
\delta.
\] Therefore $P_t^{\mathbf{E}}\mathbf{u}_j(x)\to \mathbf{K}_t^{\mathbf{E}}(x,y)\boldsymbol{\eta}$. The same argument applied to the scalar kernel gives $P_t\rho_j(x)\to p_t(x,y)$. Passing to the limit, \[
|\mathbf{K}_t^{\mathbf{E}}(x,y)\boldsymbol{\eta}|_{\mathbf{h}_{\mathbf{E}}}
\leq
p_t(x,y)|\boldsymbol{\eta}|_{\mathbf{h}_{\mathbf{E}}}.
\] Taking the supremum over unit vectors in $\mathbf{E}_y$ gives the operator bound. \end{proof}

\begin{remark}[Scope of the kernel construction] The preceding construction gives a smooth kernel for every $t>0$ and the domination estimate $\|\mathbf K_t^{\mathbf E}(x,y)\|\leq p_t(x,y)$. Estimate \eqref{eq:regularizacion-calor-cota-tiempo-compacto} describes regularization on compact sets and is uniform when time stays bounded away from zero. Asymptotic expansions as $t\to 0^{+}$ additionally require the local analysis in Chapter~\ref{cap:operadores-pseudodiferenciales}. For the scalar construction of these expansions, see \cite[Chapter VI, §§3--4]{Chavel1984Eigenvalues}. \end{remark}

The following scalar properties will be needed to extend the semigroup to the entire $L^p(M)$ scale. \begin{proposition}[Properties of the scalar kernel] \label{prop:propiedades-nucleo-calor-escalar} Let $(M,\mathbf{g})$ be a connected, complete Riemannian manifold without boundary. The scalar kernel $p_t(x,y)$ is smooth, symmetric, and nonnegative. It satisfies \[
p_{t+s}(x,y)
=
\int_Mp_t(x,z)p_s(z,y)\,d\lambda_{\mathbf{g}}(z)
\] and \[
\int_Mp_t(x,y)\,d\lambda_{\mathbf{g}}(y)\leq1.
\] \end{proposition}

\begin{proof} Smoothness, symmetry, and the composition identity are the trivial-bundle case of the preceding theorem. Specify the passage from almost-everywhere bounds to pointwise bounds. If $f\in C_c^\infty(M)$ and $f\geq0$, Corollary~\ref{cor:semigrupo-escalar-submarkoviano} gives $P_tf\geq0$ almost everywhere. Its smooth representative is continuous, and every nonempty open set has positive Riemannian measure; hence $P_tf(x)\geq0$ for every $x$. If also $f\leq1$, the same argument gives $P_tf(x)\leq1$ for every $x$. The integral representation applied to nonnegative test functions also shows that $p_t$ is real: if its imaginary part were nonzero at a point, a test function supported where that part has constant sign would give a nonzero imaginary part of $P_tf(x)$.

Suppose that $p_t(x_0,y_0)<0$. By continuity, there exist a neighborhood $W$ of $y_0$ and a constant $c>0$ such that $p_t(x_0,y)\leq-c$ for every $y\in W$. Choose a nonzero function $f\in C_c^\infty(W)$ with $f\geq0$. Then \[
P_tf(x_0)
=
\int_Wp_t(x_0,y)f(y)\,d\lambda_{\mathbf{g}}(y)
\leq
-c\int_Wf\,d\lambda_{\mathbf{g}}
<0,
\] contradicting positivity of the scalar semigroup proved in Corollary~\ref{cor:semigrupo-escalar-submarkoviano}. Therefore $p_t\geq0$.

Take the functions $(\chi_j)$ from Lemma~\ref{lem:funciones-corte-completitud-capitulo-calor}. We chose this sequence so that $0\leq\chi_j\leq1$, it is increasing, and it converges pointwise to $1$. Contractivity on $L^\infty(M)$ gives $0\leq P_t\chi_j(x)\leq1$. Since $p_t(x,y)\geq0$, the sequence $p_t(x,\cdot)\chi_j$ is also increasing and converges pointwise to $p_t(x,\cdot)$. The monotone convergence theorem~\ref{convergencia monotona} gives \[
\int_Mp_t(x,y)\,d\lambda_{\mathbf{g}}(y)
=
\lim_{j\to\infty}\int_Mp_t(x,y)\chi_j(y)\,d\lambda_{\mathbf{g}}(y)
=
\lim_{j\to\infty}P_t\chi_j(x)
\leq1.
\] \end{proof}

\section{Extension of the semigroup to \texorpdfstring{$L^p(M,\mathbf{E})$}{Lp(E)}}

\subsection{The scalar semigroup on the \texorpdfstring{$L^p(M)$}{Lp} scale}

\begin{proposition}[Scalar semigroup on $L^p(M)$] \label{prop:semigrupo-escalar-Lp-calor} Let $(M,\mathbf{g})$ be a connected, complete Riemannian manifold without boundary. For each $1\leq p\leq\infty$ and $t\geq0$, the operator $P_t$ extends to a positive contraction $P_t^p$ on $L^p(M)$. The extensions are consistent. If $1\leq p<\infty$, $(P_t^p)_{t\geq0}$ is a $C_0$-semigroup. \end{proposition}

\begin{proof} The bound on $L^\infty(M)$ was proved in Corollary~\ref{cor:semigrupo-escalar-submarkoviano}. For $f\in L^1(M)\cap L^2(M)$, self-adjointness and the bound on $L^\infty(M)$ give \[
\|P_tf\|_{L^1(M)}
=
\sup_{\substack{\varphi\in L^2(M)\cap L^\infty(M)\\ \|\varphi\|_{L^\infty(M)}\leq1}}
\left|\int_M(P_tf)\varphi\,d\lambda_{\mathbf{g}}\right|
=
\sup_{\substack{\varphi\in L^2(M)\cap L^\infty(M)\\ \|\varphi\|_{L^\infty(M)}\leq1}}
\left|\int_Mf(P_t\varphi)\,d\lambda_{\mathbf{g}}\right|
\leq
\|f\|_{L^1(M)}.
\] In the first line, the norm is initially understood to take values in $[0,\infty]$; the last bound proves a posteriori that $P_tf\in L^1(M)$. Let us justify the supremum formula without assuming this conclusion. If $h\in L^2(M)$, choose an increasing exhaustion $(K_j)_{j\in\mathbb N}$ by compact sets; each $K_j$ has finite measure. In the real case, set $\varphi_j=\mathbf 1_{K_j}\operatorname{sgn}(h)$ and, in the complex case, \[
\varphi_j(x):=
\begin{cases}
\mathbf 1_{K_j}(x)\,\frac{\overline{h(x)}}{|h(x)|},&h(x)\neq0,\\
0,&h(x)=0.
\end{cases}
\] Then $\varphi_j\in L^2(M)\cap L^\infty(M)$, $\|\varphi_j\|_{L^\infty(M)}\leq1$, and \[
\left|\int_Mh \varphi_j\,d\lambda_{\mathbf{g}}\right|
=
\int_{K_j}|h|\,d\lambda_{\mathbf{g}}
\longrightarrow
\int_M|h|\,d\lambda_{\mathbf{g}}
\] by monotone convergence, even if the value is infinite. The reverse inequality holds for every $\varphi$ by the definition of the extended integral. The Riesz--Thorin theorem~\ref{teo:riesz-thorin-interpolacion} gives contractivity on $L^p(M)$ for $1<p<\infty$. Density of $L^2(M)\cap L^p(M)$ in $L^p(M)$ when $p<\infty$ gives the unique extension. For $p=\infty$, the integral formula with $p_t\geq0$ and mass at most one defines the extension on $L^\infty(M)$.

The integral formula also represents the extensions for finite exponents. Indeed, if $f\in L^p(M)$ and $1\leq p<\infty$, Jensen's inequality for the subprobability measure $p_t(x,y)\,d\lambda_{\mathbf g}(y)$ gives \[
 \left(\int_Mp_t(x,y)|f(y)|\,d\lambda_{\mathbf g}(y)\right)^p
 \leq\int_Mp_t(x,y)|f(y)|^p\,d\lambda_{\mathbf g}(y).
\] For $p=1$, this inequality is an equality. Tonelli's theorem, symmetry of the kernel, and its mass bound show that the integral of the right-hand side with respect to $x$ does not exceed $\|f\|_{L^{p}(M)}^p$. Thus the kernel integral converges absolutely for almost every $x$ and defines a contraction on $L^p(M)$. It agrees with $P_t$ on the intersection with $L^2(M)$, and density identifies this contraction with $P_t^p$.

For consistency, let $p,q<\infty$ and $f\in L^p(M)\cap L^q(M)$. Using the preceding exhaustion $(K_j)$, set \[
 f_j=\mathbf1_{K_j}\mathbf1_{\{|f|\leq j\}}f.
\] Each $f_j$ is bounded and has support of finite measure; hence it belongs simultaneously to $L^2(M)$, $L^p(M)$, and $L^q(M)$. Dominated convergence applied to $|f_j-f|^p$ and $|f_j-f|^q$ gives convergence in both spaces. By contractivity, $P_t^pf_j\to P_t^pf$ in $L^p(M)$ and $P_t^qf_j\to P_t^qf$ in $L^q(M)$. The terms of the two sequences agree for every $j$ because they arise from the same operator on $L^2(M)$. Pairing with a compactly supported smooth test function, Hölder's inequality identifies their distributional limits; thus $P_t^pf=P_t^qf$ almost everywhere. If one exponent is infinite, the integral formula already proved for each exponent directly gives the same function.

For $p<\infty$, now take $f_j\to f$ in $L^p(M)$ with $f_j\in L^2(M)\cap L^p(M)$. Consistency allows us to apply the semigroup identity on $L^2(M)$ to each $f_j$, and \[
 \|P_t^pP_s^pf-P_{t+s}^pf\|_{L^p(M)}
 \leq\|P_t^pP_s^p(f-f_j)\|_{L^p(M)}
      +\|P_{t+s}^p(f_j-f)\|_{L^p(M)}
 \leq2\|f-f_j\|_{L^p(M)}\longrightarrow0.
\] For $p=\infty$ and $s,t>0$, the double integral arising in the composition of kernels converges absolutely, since \[
 \int_M\int_Mp_t(x,z)p_s(z,y)|f(y)|
      \,d\lambda_{\mathbf g}(y)\,d\lambda_{\mathbf g}(z)
 \leq\|f\|_{L^\infty(M)}.
\] Fubini's theorem and the Chapman--Kolmogorov identity for the kernel give $P_t^\infty P_s^\infty f=P_{t+s}^\infty f$. When one time is zero, the equality follows from $P_0^p=I$.

Let us prove strong continuity for $p<\infty$. If $f\in C_c^\infty(M)$, then $Lf\in L^p(M)\cap L^2(M)$. In $L^2(M)$, spectral differentiability gives \[
P_t f-f
=
-\int_0^tP_s(Lf)\,ds.
\] Let us justify strong measurability of the integrand before using the strong continuity we are proving. Set $g=Lf$ and use the cutoffs $\chi_j$ from Lemma~\ref{lem:funciones-corte-completitud-capitulo-calor}. For $s>0$, the map $s\mapsto\chi_jP_sg$ is continuous with values in $L^p(M)$: heat regularization gives uniform continuity on the compact support of $\chi_j$ on each compact time interval in $(0,\infty)$. For every $s>0$, dominated convergence in $M$ gives $\chi_jP_sg\to P_s^pg$ in $L^p(M)$. Thus $s\mapsto P_s^pg$ is the pointwise limit of strongly measurable maps; assigning its value at $s=0$ does not change this property. Moreover, its norm does not exceed $\|g\|_{L^p(M)}$. The Bochner integral exists in $L^p(M)$. Pairing the integrals in $L^2(M)$ and $L^p(M)$ with compactly supported smooth test functions, consistency gives the same result; their distributional realizations agree. Therefore the same identity holds in $L^p(M)$: \[
P_t^pf-f
=
-\int_0^tP_s^p(Lf)\,ds.
\] Thus \[
\|P_t^pf-f\|_{L^p(M)}
\leq
\int_0^t\|P_s^p(Lf)\|_{L^p(M)}\,ds
\leq
t\|Lf\|_{L^p(M)}
\longrightarrow0
\] as $t\to0^+$. For $f\in L^p(M)$, take $g_j\in C_c^\infty(M)$ with $g_j\to f$ in $L^p(M)$. Then \[
 \|P_t^pf-f\|_{L^p(M)}
 \leq2\|f-g_j\|_{L^p(M)}+t\|Lg_j\|_{L^p(M)}.
\] First fix $j$ to control the first term, and then let $t\to0^+$. This proves strong continuity at zero for each $f$. The semigroup property and contractivity transfer it to any time: for $h>0$, both the right and left differences at a positive time are bounded by $\|P_h^pf-f\|_{L^p(M)}$. \end{proof}

\subsection{The covariant semigroup on \texorpdfstring{$L^p(M,\mathbf{E})$}{Lp(E)}}

\begin{theorem}[Covariant semigroup on $L^p(M,\mathbf E)$] \label{teo:semigrupo-calor-covariante-Lp} Let $(M,\mathbf{g})$ be a connected, complete Riemannian manifold without boundary, and let $\mathbf{E}\to M$ be a smooth vector bundle of finite rank with a bundle metric and a compatible connection. Let $1\leq p\leq\infty$. For each $t\geq0$, $P_t^{\mathbf{E}}$ admits a contractive extension \[
P_t^{\mathbf{E},p}\colon L^p(M,\mathbf{E})\longrightarrow L^p(M,\mathbf{E}).
\] The extensions are consistent and satisfy: \begin{enumerate}[label=(\alph*)] \item $P_{t+s}^{\mathbf{E},p}=P_t^{\mathbf{E},p}P_s^{\mathbf{E},p}$; \item if $1\leq p<\infty$, $(P_t^{\mathbf{E},p})_{t\geq0}$ is a $C_0$-semigroup; \item if $t>0$ and $\mathbf{u}\in L^p(M,\mathbf{E})$, \[
P_t^{\mathbf{E},p}\mathbf{u}(x)
=
\int_M\mathbf{K}_t^{\mathbf{E}}(x,y)\mathbf{u}(y)\,d\lambda_{\mathbf{g}}(y)
\] for almost every $x$; \item for $t\geq0$, \[
|P_t^{\mathbf{E},p}\mathbf{u}|_{\mathbf{h}_{\mathbf{E}}}
\leq
P_t^p(|\mathbf{u}|_{\mathbf{h}_{\mathbf{E}}})
\] almost everywhere. \end{enumerate} For $p=\infty$, the preceding conclusions do not include strong continuity on $t=0$. \end{theorem}

\begin{proof} If $\mathbf{u}\in L^2(M,\mathbf{E})\cap L^p(M,\mathbf{E})$, the Kato--Simon inequality and scalar contractivity give \[
\|P_t^{\mathbf{E}}\mathbf{u}\|_{L^p(M,\mathbf{E})}
\leq
\|P_t^p(|\mathbf{u}|_{\mathbf{h}_{\mathbf{E}}})\|_{L^p(M)}
\leq
\|\mathbf{u}\|_{L^p(M,\mathbf{E})}.
\] For $p<\infty$, the subspace $L^2(M,\mathbf{E})\cap L^p(M,\mathbf{E})$ is dense in $L^p(M,\mathbf{E})$, and Theorem~\ref{teo:extension-operadores-subespacio-denso} provides the unique extension. The pointwise inequality passes to the limit: if $\mathbf{u}_j\to \mathbf{u}$ in $L^p(M,\mathbf{E})$, then, after passing to a subsequence, $P_t^{\mathbf{E},p}\mathbf{u}_j\to P_t^{\mathbf{E},p}\mathbf{u}$ and $P_t^p(|\mathbf{u}_j|_{\mathbf{h}_{\mathbf{E}}})\to P_t^p(|\mathbf{u}|_{\mathbf{h}_{\mathbf{E}}})$ almost everywhere.

For $p=\infty$ and $t>0$, define \[
P_t^{\mathbf{E},\infty}\mathbf{u}(x)
:=
\int_M\mathbf{K}_t^{\mathbf{E}}(x,y)\mathbf{u}(y)\,d\lambda_{\mathbf{g}}(y),
\] and set $P_0^{\mathbf{E},\infty}:=I$. The kernel bound and Proposition~\ref{prop:propiedades-nucleo-calor-escalar} give \[
|P_t^{\mathbf{E},\infty}\mathbf{u}(x)|_{\mathbf{h}_{\mathbf{E}}}
\leq
\|\mathbf{u}\|_{L^\infty(M,\mathbf{E})}
\int_Mp_t(x,y)\,d\lambda_{\mathbf{g}}(y)
\leq
\|\mathbf{u}\|_{L^\infty(M,\mathbf{E})}.
\] The same estimate proves absolute convergence of the integral and pointwise domination.

For $p<\infty$ and $\mathbf u\in L^p(M,\mathbf E)$, kernel domination and the Jensen--Tonelli argument in the scalar proposition show that \[
 \int_M\|\mathbf K_t^{\mathbf E}(x,y)\mathbf u(y)\|_{\mathbf h_{\mathbf E}}
       \,d\lambda_{\mathbf g}(y)
 \leq P_t^p(|\mathbf u|_{\mathbf h_{\mathbf E}})(x)<\infty
\] for almost every $x$. This integral defines a contractive operator on $L^p(M,\mathbf E)$ agreeing with $P_t^{\mathbf E}$ on $L^2(M,\mathbf E)\cap L^p(M,\mathbf E)$. The approximation is obtained by taking \[
 \mathbf u_j
 =\mathbf1_{K_j}\mathbf1_{\{|\mathbf u|_{\mathbf h_{\mathbf E}}\leq j\}}
       \mathbf u.
\] Then $\mathbf u_j\in L^2(M,\mathbf E)\cap L^p(M,\mathbf E)$ and $\|\mathbf u_j-\mathbf u\|_{L^p(M,\mathbf E)}\to0$. Both the integral operator and $P_t^{\mathbf E,p}$ send this difference to a section whose norm does not exceed $\|\mathbf u_j-\mathbf u\|_{L^p(M,\mathbf E)}$. Their limits agree, proving the integral formula for every $\mathbf u$.

Since all exponents are represented by the same kernel, the extensions are consistent, including the infinite exponent. For $s,t>0$, the norm of the integrand in the composition of operators is dominated by \[
 p_t(x,z)p_s(z,y)|\mathbf u(y)|_{\mathbf h_{\mathbf E}}.
\] The double integral of this function is $P_{t+s}^p(|\mathbf u|_{\mathbf h_{\mathbf E}})(x)$, by Tonelli's theorem and the scalar Chapman--Kolmogorov identity. It is finite for almost every $x$ when $p<\infty$, and for every $x$ when $p=\infty$. Fubini's theorem in the finite-dimensional fiber $\mathbf E_x$ therefore allows the integrals to be interchanged. The composition identity for the covariant kernel gives $P_t^{\mathbf E,p}P_s^{\mathbf E,p}\mathbf u
=P_{t+s}^{\mathbf E,p}\mathbf u$ almost everywhere. Zero times are included by setting $P_0^{\mathbf E,p}=I$.

Finally, if $\mathbf{u}\in\Gamma_c(\mathbf{E})$, then $L_{\mathbf{E},0}\mathbf{u}\in L^p(M,\mathbf{E})\cap L^2(M,\mathbf{E})$. The identity \[
P_t^{\mathbf{E}}\mathbf{u}-\mathbf{u}
=
-\int_0^tP_s^{\mathbf{E}}(L_{\mathbf{E},0}\mathbf{u})\,ds
\] is obtained first in $L^2(M,\mathbf{E})$. Strong measurability in $L^p(M,\mathbf E)$ is checked as in Proposition~\ref{prop:semigrupo-escalar-Lp-calor}: for each compact cutoff $\chi_j$, covariant regularization makes $s\mapsto\chi_jP_s^{\mathbf E}(L_{\mathbf E,0}\mathbf u)$ continuous when $s>0$, and these maps converge pointwise in $L^p(M,\mathbf E)$ to the integrand. Contractivity bounds its norm by $\|L_{\mathbf E,0}\mathbf u\|_{L^p(M,\mathbf E)}$, so it is Bochner integrable. Consistency, checked against test sections, identifies this integral with the one obtained in $L^2(M,\mathbf E)$. Therefore \[
P_t^{\mathbf{E},p}\mathbf{u}-\mathbf{u}
=
-\int_0^tP_s^{\mathbf{E},p}(L_{\mathbf{E},0}\mathbf{u})\,ds.
\] Thus \[
\|P_t^{\mathbf{E},p}\mathbf{u}-\mathbf{u}\|_{L^p(M,\mathbf{E})}
\leq
t\|L_{\mathbf{E},0}\mathbf{u}\|_{L^p(M,\mathbf{E})}
\longrightarrow0
\] as $t\to0^+$. For an arbitrary section $\mathbf u\in L^p(M,\mathbf E)$, choose $\mathbf v_j\in\Gamma_c(\mathbf E)$ with $\mathbf v_j\to\mathbf u$ in $L^p(M,\mathbf E)$ and use \[
 \|P_t^{\mathbf E,p}\mathbf u-\mathbf u\|_{L^p(M,\mathbf E)}
 \leq2\|\mathbf u-\mathbf v_j\|_{L^p(M,\mathbf E)}
           +t\|L_{\mathbf E,0}\mathbf v_j\|_{L^p(M,\mathbf E)}.
\] First fix $j$ and then let $t\to0^+$; finally, let $j\to\infty$. The semigroup property transfers strong continuity to all times, just as in the scalar case. \end{proof}

\begin{definition}[Generator on $L^p(M,\mathbf E)$] \label{def:generador-laplaciano-conexion-Lp} Let $(M,\mathbf{g})$ be a connected, complete Riemannian manifold without boundary, and let $\mathbf{E}\to M$ be a smooth vector bundle of finite rank with a bundle metric and a compatible connection. \index{infinitesimal generator} Let $1\leq p<\infty$. Define \[
\mathcal D(L_{\mathbf{E},p})
:=
\left\{\mathbf{u}\in L^p(M,\mathbf{E})\middle|
\lim_{t\to0^+}\frac{\mathbf{u}-P_t^{\mathbf{E},p}\mathbf{u}}{t}
\text{ exists in }L^p(M,\mathbf{E})
\right\}
\] and \[
L_{\mathbf{E},p}\mathbf{u}
:=
\lim_{t\to0^+}\frac{\mathbf{u}-P_t^{\mathbf{E},p}\mathbf{u}}{t}.
\] Then $-L_{\mathbf{E},p}$ is the generator of $(P_t^{\mathbf{E},p})_{t\geq0}$. \end{definition}

\begin{proposition} \label{prop:generador-Lp-extiende-laplaciano} Let $(M,\mathbf{g})$ be a connected, complete Riemannian manifold without boundary, and let $\mathbf{E}\to M$ be a smooth vector bundle of finite rank with a bundle metric and a compatible connection. For $1\leq p<\infty$, $L_{\mathbf{E},p}$ is closed, densely defined, and maximally accretive. Moreover, \[
\Gamma_c(\mathbf{E})\subseteq\mathcal D(L_{\mathbf{E},p}),
\qquad
L_{\mathbf{E},p}\mathbf{u}=L_{\mathbf{E},0}\mathbf{u}
\quad\text{for }\mathbf{u}\in\Gamma_c(\mathbf{E}),
\] and $L_{\mathbf{E},2}=L_{\mathbf{E}}$. \end{proposition}

\begin{proof} Closedness and density of the domain follow from Proposition~\ref{prop:propiedades-generador-semigrupo-C0}; maximal accretivity is the Lumer--Phillips theorem, Theorem~\ref{teo:lumer-phillips}, applied to the generator $-L_{\mathbf{E},p}$. If $\mathbf{u}\in\Gamma_c(\mathbf{E})$, then \[
\frac{\mathbf{u}-P_t^{\mathbf{E},p}\mathbf{u}}{t}
=
\frac{1}{t}\int_0^tP_s^{\mathbf{E},p}(L_{\mathbf{E},0}\mathbf{u})\,ds.
\] Strong continuity of the semigroup implies that the right-hand side converges to $L_{\mathbf{E},0}\mathbf{u}$ in $L^p(M,\mathbf{E})$ as $t\to0^+$. For $p=2$, the semigroups defined by extension and by spectral calculus agree; uniqueness of the generator gives $L_{\mathbf{E},2}=L_{\mathbf{E}}$. \end{proof}

\section{Bessel potentials and spaces of real order}

\subsection{Bessel potentials}

The functional calculus of the Laplacian transfers the Euclidean multiplier $(1+\|\xi\|^2)^{-\displaystyle\frac{s}{2}}$ to the manifold. The resulting Bessel potentials encode a gain of regularity without choosing coordinates.

\begin{definition}[Bessel potential] \label{def:sobolev-fraccionario-y-nucleo-calor-el-laplaciano-en-haces-vectoriales-y-el-nucleo-d} Let $(M,\mathbf{g})$ be a connected, complete Riemannian manifold without boundary, and let $\mathbf{E}\to M$ be a smooth vector bundle of finite rank with a bundle metric and a compatible connection. \index{Bessel potential} Let $s>0$ and $1\leq p<\infty$. Define \[
J_{s,p}^{\mathbf{E}}\mathbf{u}
:=
\frac{1}{\Gamma(\frac{s}{2})}
\int_0^\infty
t^{\frac{s}{2}-1}e^{-t}P_t^{\mathbf{E},p}\mathbf{u}\,dt,
\qquad \mathbf{u}\in L^p(M,\mathbf{E}),
\] where the integral is a Bochner integral in $L^p(M,\mathbf{E})$. \end{definition}

\begin{proposition}[Properties of Bessel potentials] \label{prop:propiedades-potenciales-bessel-haces} Let $(M,\mathbf{g})$ be a connected, complete Riemannian manifold without boundary, and let $\mathbf{E}\to M$ be a smooth vector bundle of finite rank with a bundle metric and a compatible connection. Let $r,s>0$ and $1\leq p,q<\infty$. \begin{enumerate}[label=(\alph*)] \item $J_{s,p}^{\mathbf{E}}\in\mathcal L(L^p(M,\mathbf{E}))$ and $\|J_{s,p}^{\mathbf{E}}\|_{\mathcal L(L^p(M,\mathbf{E}))}\leq1$; \item $J_{s,p}^{\mathbf{E}}J_{r,p}^{\mathbf{E}}=J_{s+r,p}^{\mathbf{E}}$; \item $J_{2,p}^{\mathbf{E}}=(I+L_{\mathbf{E},p})^{-1}$ and \[
J_{2m,p}^{\mathbf{E}}=(I+L_{\mathbf{E},p})^{-m}
\] for $m\in\mathbb N$; \item $J_{s,p}^{\mathbf{E}}$ is injective and its range is dense in $L^p(M,\mathbf{E})$; \item if $\mathbf{u}\in L^p(M,\mathbf{E})\cap L^q(M,\mathbf{E})$, then $J_{s,p}^{\mathbf{E}}\mathbf{u}=J_{s,q}^{\mathbf{E}}\mathbf{u}$ almost everywhere. \end{enumerate} \end{proposition}

\begin{proof} Contractivity of the semigroup gives \[
\int_0^\infty
\left\|
\frac{t^{\frac{s}{2}-1}e^{-t}}{\Gamma(\frac{s}{2})}P_t^{\mathbf{E},p}\mathbf{u}
\right\|_{L^p(M,\mathbf{E})}dt
\leq
\frac{\|\mathbf{u}\|_{L^p(M,\mathbf{E})}}{\Gamma(\frac{s}{2})}
\int_0^\infty t^{\frac{s}{2}-1}e^{-t}\,dt
=
\|\mathbf{u}\|_{L^p(M,\mathbf{E})}.
\] This proves convergence of the Bochner integral and $\|J_{s,p}^{\mathbf{E}}\|_{\mathcal L(L^p(M,\mathbf{E}))}\leq1$.

For the semigroup law, the norm of the double integrand is bounded by \[
\frac{t^{\frac{s}{2}-1}\tau^{\frac{r}{2}-1}e^{-(t+\tau)}}{\Gamma(\frac{s}{2})\Gamma(\frac{r}{2})}
\|\mathbf{u}\|_{L^p(M,\mathbf{E})}.
\] The integral of this function over $(0,\infty)^2$ is $\|\mathbf{u}\|_{L^p(M,\mathbf{E})}$. Thus Fubini's theorem for Bochner integrals~\ref{teo:b5-fubini-bochner} applies and \[
J_{s,p}^{\mathbf{E}}J_{r,p}^{\mathbf{E}}\mathbf{u}
=
\frac{1}{\Gamma(\frac{s}{2})\Gamma(\frac{r}{2})}
\int_0^\infty\int_0^\infty
t^{\frac{s}{2}-1}\tau^{\frac{r}{2}-1}e^{-(t+\tau)}P_{t+\tau}^{\mathbf{E},p}\mathbf{u}
\,d\tau\,dt.
\] Make the change of variables \[
a=t+\tau,
\qquad
\theta=\frac{t}{t+\tau}.
\] Its Jacobian is $a$, and the beta--gamma identity in Lemma~\ref{lem:identidad-beta-gamma} gives \[
\int_0^1\theta^{\frac{s}{2}-1}(1-\theta)^{\frac{r}{2}-1}\,d\theta
=
\frac{\Gamma(\frac{s}{2})\Gamma(\frac{r}{2})}{\Gamma\left(\frac{s+r}{2}\right)}.
\] We obtain $J_{s+r,p}^{\mathbf{E}}\mathbf{u}$.

The resolvent formula, Proposition~\ref{prop:resolvente-transformada-laplace-semigrupo}, applied to the generator $-L_{\mathbf{E},p}$ with $\lambda=1$, gives $J_{2,p}^{\mathbf{E}}=(I+L_{\mathbf{E},p})^{-1}$. The semigroup law gives the formula for $J_{2m,p}^{\mathbf{E}}$.

If $J_{s,p}^{\mathbf{E}}\mathbf{u}=0$, choose $m\in\mathbb N$ with $2m>s$. Then \[
J_{2m,p}^{\mathbf{E}}\mathbf{u}
=
J_{2m-s,p}^{\mathbf{E}}J_{s,p}^{\mathbf{E}}\mathbf{u}
=0.
\] The operator $(I+L_{\mathbf{E},p})^{-m}$ is injective, so $\mathbf{u}=0$. For density, observe that \[
J_{2m,p}^{\mathbf{E}}
=
J_{s,p}^{\mathbf{E}}J_{2m-s,p}^{\mathbf{E}}.
\] Therefore \[
\operatorname{Ran}(J_{2m,p}^{\mathbf{E}})
\subseteq
\operatorname{Ran}(J_{s,p}^{\mathbf{E}}).
\] Moreover, \[
\operatorname{Ran}(J_{2m,p}^{\mathbf{E}})
=
\operatorname{Ran}\bigl((I+L_{\mathbf{E},p})^{-m}\bigr)
=
\mathcal D\bigl((I+L_{\mathbf{E},p})^m\bigr),
\] where the last space is viewed with its inclusion into $L^p(M,\mathbf{E})$. The domain of a power of the generator is dense by Proposition~\ref{prop:propiedades-generador-semigrupo-C0}. Consequently, $\operatorname{Ran}(J_{s,p}^{\mathbf{E}})$ is dense in $L^p(M,\mathbf{E})$.

Finally, consistency of the semigroups implies that, for $\mathbf{u}\in L^p(M,\mathbf{E})\cap L^q(M,\mathbf{E})$, the truncated integrals \[
J_{s}^{\mathbf{E},[a,b]}\mathbf{u}
:=
\frac{1}{\Gamma(\frac{s}{2})}
\int_a^b t^{\frac{s}{2}-1}e^{-t}P_t^{\mathbf{E}}\mathbf{u}\,dt
\] are the same sections for all $0<a<b<\infty$. As $a\to0^+$ and $b\to+\infty$, these sections converge to $J_{s,p}^{\mathbf{E}}\mathbf{u}$ in $L^p(M,\mathbf{E})$ and to $J_{s,q}^{\mathbf{E}}\mathbf{u}$ in $L^q(M,\mathbf{E})$. From any sequence of intervals exhausting $(0,+\infty)$, one can extract a subsequence converging almost everywhere in both spaces: first extract one for convergence in $L^p(M,\mathbf{E})$, and then a further subsequence for convergence in $L^q(M,\mathbf{E})$. The two pointwise limits agree, and hence $J_{s,p}^{\mathbf{E}}\mathbf{u}=J_{s,q}^{\mathbf{E}}\mathbf{u}$ almost everywhere. \end{proof}

\subsection{The Bessel kernel}

The spectral representation of the potential can be converted into an integral representation. Its kernel separates the singularity near the diagonal from the decay at large distances and makes the regularization properties visible.

\begin{definition}[Covariant Bessel kernel] \label{def:nucleo-bessel-covariante} Let $(M,\mathbf{g})$ be a connected, complete Riemannian manifold without boundary, and let $\mathbf{E}\to M$ be a smooth vector bundle of finite rank with a bundle metric and a compatible connection. For $s>0$, define \[
G_s^{\mathbf{E}}(x,y)
:=
\frac{1}{\Gamma(\frac{s}{2})}
\int_0^\infty t^{\frac{s}{2}-1}e^{-t}\mathbf{K}_t^{\mathbf{E}}(x,y)\,dt
\] whenever the integral converges. Also define \[
g_s(x,y)
:=
\frac{1}{\Gamma(\frac{s}{2})}
\int_0^\infty t^{\frac{s}{2}-1}e^{-t}p_t(x,y)\,dt,
\] initially allowing the value $+\infty$. \end{definition}

\begin{proposition} \label{prop:representacion-nucleo-bessel-covariante} Let $(M,\mathbf{g})$ be a connected, complete Riemannian manifold without boundary, and let $\mathbf{E}\to M$ be a smooth vector bundle of finite rank with a bundle metric and a compatible connection. For almost every $(x,y)\in M\times M$, $g_s(x,y)<\infty$, the integral defining $G_s^{\mathbf{E}}(x,y)$ converges and \[
\|G_s^{\mathbf{E}}(x,y)\|_{\mathcal L((\mathbf{E}_y,\mathbf{h}_{\mathbf{E}}(y)),(\mathbf{E}_x,\mathbf{h}_{\mathbf{E}}(x)))}
\leq
g_s(x,y).
\] Moreover, \[
\int_Mg_s(x,y)\,d\lambda_{\mathbf{g}}(y)\leq1
\] for every $x$, and, if $\mathbf{u}\in L^p(M,\mathbf{E})$ with $1\leq p<\infty$, \[
J_{s,p}^{\mathbf{E}}\mathbf{u}(x)
=
\int_MG_s^{\mathbf{E}}(x,y)\mathbf{u}(y)\,d\lambda_{\mathbf{g}}(y)
\] for almost every $x$. \end{proposition}

\begin{proof} The function $(t,x,y)\mapsto p_t(x,y)$ is nonnegative and measurable. Thus $g_s$ is measurable as the integral of a nonnegative function. For each $x\in M$, Tonelli's theorem~\ref{teo:tonelli} and the mass bound for $p_t$ give \begin{align*}
\int_Mg_s(x,y)\,d\lambda_{\mathbf{g}}(y)
&=
\frac{1}{\Gamma(\frac{s}{2})}
\int_0^\infty t^{\frac{s}{2}-1}e^{-t}
\left(\int_Mp_t(x,y)\,d\lambda_{\mathbf{g}}(y)\right)dt\\
&\leq
\frac{1}{\Gamma(\frac{s}{2})}
\int_0^\infty t^{\frac{s}{2}-1}e^{-t}\,dt
=1.
\end{align*} Since $p_t(x,y)=p_t(y,x)$, also \begin{equation}
\label{eq:simetria-nucleo-bessel-escalar}
g_s(x,y)=g_s(y,x).
\end{equation} In particular, $g_s(x,\cdot)$ is finite almost everywhere for each $x$. Since the set \[
N_s:=\{(x,y)\in M\times M\mid g_s(x,y)=+\infty\}
\] is measurable, Tonelli's theorem applied to $\mathbf 1_{N_s}$ shows that $(\lambda_{\mathbf{g}}\otimes\lambda_{\mathbf{g}})(N_s)=0$.

If $(x,y)\notin N_s$, the bound \[
\|\mathbf{K}_t^{\mathbf{E}}(x,y)\|_{\mathcal L((\mathbf{E}_y,\mathbf{h}_{\mathbf{E}}(y)),(\mathbf{E}_x,\mathbf{h}_{\mathbf{E}}(x)))}
\leq p_t(x,y)
\] proves that the function with values in the finite-dimensional space $\operatorname{Hom}(\mathbf{E}_y,\mathbf{E}_x)$ appearing in the definition of $G_s^{\mathbf{E}}(x,y)$ is absolutely integrable. Thus the integral converges and \[
\|G_s^{\mathbf{E}}(x,y)\|_{\mathcal L((\mathbf{E}_y,\mathbf{h}_{\mathbf{E}}(y)),(\mathbf{E}_x,\mathbf{h}_{\mathbf{E}}(x)))}
\leq g_s(x,y).
\] On $N_s$, we henceforth assign the value $G_s^{\mathbf{E}}(x,y)=0$. This choice gives a measurable kernel. Indeed, on the product of two trivializing open subsets, the components of \[
\frac{1}{\Gamma(\frac{s}{2})}
\int_{\frac{1}{\nu}}^{\nu}
 t^{\frac{s}{2}-1}e^{-t}\mathbf{K}_t^{\mathbf{E}}(x,y)\,dt,
\qquad \nu\in\mathbb N,
\] are measurable, in fact smooth, functions. Outside $N_s$, they converge componentwise to those of $G_s^{\mathbf{E}}$, while $N_s$ is measurable. Thus the limit defined on $M\times M\setminus N_s$ and extended by zero to $N_s$ is measurable in every pair of local frames. This proves the first assertions.

Before treating the covariant kernel, we record boundedness of the scalar integral operator. For a measurable function $f\geq0$, set \[
(\mathcal T_sf)(x):=\int_Mg_s(x,y)f(y)\,d\lambda_{\mathbf{g}}(y),
\] initially allowing the value $+\infty$. If $1\leq p<\infty$, Jensen's inequality for the subprobability measure $g_s(x,y)\,d\lambda_{\mathbf{g}}(y)$ gives \[
|\mathcal T_sf(x)|^p
\leq
\int_Mg_s(x,y)|f(y)|^p\,d\lambda_{\mathbf{g}}(y).
\] For $p=1$, the same inequality is immediate. Integrate in $x$ and use Tonelli's theorem, the symmetry \eqref{eq:simetria-nucleo-bessel-escalar}, and the mass bound already proved. We obtain \begin{equation}
\label{eq:contractividad-operador-integral-bessel-escalar}
\|\mathcal T_sf\|_{L^p(M)}
\leq
\|f\|_{L^p(M)}.
\end{equation} Decomposition into positive and negative parts gives the same conclusion for arbitrary scalar functions. In particular, the integral defining $\mathcal T_sf$ is finite for almost every $x$ when $f\in L^p(M)$. Now identify the pointwise integral with the Bochner integral defining the potential. Put
\[
w_s(t):=\frac{t^{\frac{s}{2}-1}e^{-t}}{\Gamma(\frac{s}{2})},\qquad
\int_0^\infty w_s(t)\,dt=1.
\]
First let $\mathbf u\in\Gamma_c(\mathbf E)$. Kernel domination and Tonelli give
\[
\int_0^\infty\int_M w_s(t)
 |\mathbf K_t^{\mathbf E}(x,y)\mathbf u(y)|_{\mathbf h_{\mathbf E}}
 \,d\lambda_{\mathbf g}(y)\,dt
\leq\mathcal T_s(|\mathbf u|_{\mathbf h_{\mathbf E}})(x).
\]
The right-hand side is finite for almost every $x$ and belongs to $L^p(M)$ by \eqref{eq:contractividad-operador-integral-bessel-escalar}. At those points the fiber integral exists:
\[
\mathbf H_s(x):=\int_0^\infty w_s(t)P_t^{\mathbf E}\mathbf u(x)\,dt,
\qquad
|\mathbf H_s(x)|_{\mathbf h_{\mathbf E}}
\leq\mathcal T_s(|\mathbf u|_{\mathbf h_{\mathbf E}})(x).
\]
Set $\mathbf H_s=0$ on the remaining null set. The resulting section is measurable, as is checked in local frames using truncated integrals, and belongs to $L^p(M,\mathbf E)$.

For $\boldsymbol\psi\in\Gamma_c(\mathbf E^*)$, let $p'$ be the conjugate exponent of $p$. Hölder and semigroup contractivity give
\[
\int_0^\infty w_s(t)\int_M
 |\boldsymbol\psi(x)(P_t^{\mathbf E}\mathbf u(x))|
 \,d\lambda_{\mathbf g}(x)\,dt
\leq\|\boldsymbol\psi\|_{L^{p'}(M,\mathbf E^*)}
       \|\mathbf u\|_{L^p(M,\mathbf E)}.
\]
Pairing with $\boldsymbol\psi$ is a continuous functional on $L^p(M,\mathbf E)$. It can therefore be interchanged with the Bochner integral defining $J_{s,p}^{\mathbf E}$. Fubini, justified by the preceding bound, gives
\[
\begin{split}
\int_M\boldsymbol\psi(J_{s,p}^{\mathbf E}\mathbf u)\,d\lambda_{\mathbf g}
&=\int_0^\infty w_s(t)\int_M
       \boldsymbol\psi(P_t^{\mathbf E}\mathbf u)\,d\lambda_{\mathbf g}\,dt\\
&=\int_M\boldsymbol\psi(\mathbf H_s)\,d\lambda_{\mathbf g}.
\end{split}
\]
Uniqueness of sections representing a regular distribution implies $\mathbf H_s=J_{s,p}^{\mathbf E}\mathbf u$ almost everywhere. The absolute integrability bound in $(t,y)$ now allows Theorem~\ref{teo:b5-fubini-bochner} to be applied in each fiber, yielding
\[
\begin{split}
J_{s,p}^{\mathbf E}\mathbf u(x)
&=\int_0^\infty\int_M w_s(t)\mathbf K_t^{\mathbf E}(x,y)
               \mathbf u(y)\,d\lambda_{\mathbf g}(y)\,dt\\
&=\int_MG_s^{\mathbf E}(x,y)\mathbf u(y)\,d\lambda_{\mathbf g}(y)
\end{split}
\]
for almost every $x$. In the trivial bundle, the same argument proves $\mathcal T_sf=J_{s,p}f$ for $f\in C_c^\infty(M)$. Density and \eqref{eq:contractividad-operador-integral-bessel-escalar} extend this identity to all of $L^p(M)$.

Whenever the integral converges absolutely, define \[
\mathcal G_s^{\mathbf{E}}\mathbf{u}(x):=
\int_MG_s^{\mathbf{E}}(x,y)\mathbf{u}(y)\,d\lambda_{\mathbf{g}}(y).
\] The preceding bound gives \begin{equation}
\label{eq:dominacion-nucleo-bessel-covariante}
|\mathcal G_s^{\mathbf{E}}\mathbf{u}(x)|_{\mathbf{h}_{\mathbf{E}}}
\leq
\int_Mg_s(x,y)|\mathbf{u}(y)|_{\mathbf{h}_{\mathbf{E}}}\,d\lambda_{\mathbf{g}}(y)
=J_{s,p}(|\mathbf{u}|_{\mathbf{h}_{\mathbf{E}}})(x)
\end{equation} for almost every $x$. The right-hand side is finite almost everywhere, and its $L^p(M)$
norm does not exceed $\|\mathbf{u}\|_{L^p(M,\mathbf{E})}$. Thus $\mathcal G_s^{\mathbf{E}}$ defines a contractive operator on the subspace $\Gamma_c(\mathbf{E})$.

Now take $\mathbf{u}\in L^p(M,\mathbf{E})$ and choose $\mathbf{u}_j\in\Gamma_c(\mathbf{E})$ with $\mathbf{u}_j\to \mathbf{u}$ in $L^p(M,\mathbf{E})$. Inequality \eqref{eq:dominacion-nucleo-bessel-covariante}, applied to $\mathbf{u}_j-\mathbf{u}_k$, shows that $(\mathcal G_s^{\mathbf{E}}\mathbf{u}_j)_j$ is Cauchy in $L^p(M,\mathbf{E})$. Therefore $\mathcal G_s^{\mathbf{E}}$ extends uniquely to a contractive operator on $L^p(M,\mathbf{E})$. Since $\mathcal G_s^{\mathbf{E}}\mathbf{u}_j=J_{s,p}^{\mathbf{E}}\mathbf{u}_j$ and both operators are continuous, we obtain \[
\mathcal G_s^{\mathbf{E}}\mathbf{u}=J_{s,p}^{\mathbf{E}}\mathbf{u}
\] in $L^p(M,\mathbf{E})$.

Finally, the preceding scalar identity, now applied to $|\mathbf{u}|_{\mathbf{h}_{\mathbf{E}}}$, shows that \[
\bigl(\mathcal T_s(|\mathbf{u}|_{\mathbf{h}_{\mathbf{E}}})\bigr)(x)
=
\int_Mg_s(x,y)|\mathbf{u}(y)|_{\mathbf{h}_{\mathbf{E}}}\,d\lambda_{\mathbf{g}}(y)<\infty
\] for almost every $x$. At these points, the integral \[
\int_MG_s^{\mathbf{E}}(x,y)\mathbf{u}(y)\,d\lambda_{\mathbf{g}}(y)
\] converges absolutely. Since $\mathcal T_s$ is contractive on $L^p(M)$, $\mathcal T_s(|\mathbf{u}_j-\mathbf{u}|_{\mathbf{h}_{\mathbf{E}}})\to0$ in $L^p(M)$. Choose a subsequence, without relabeling, for which \[
\mathcal T_s(|\mathbf{u}_j-\mathbf{u}|_{\mathbf{h}_{\mathbf{E}}})(x)\longrightarrow0
\quad\text{and}\quad
\mathcal G_s^{\mathbf{E}}\mathbf{u}_j(x)\longrightarrow\mathcal G_s^{\mathbf{E}}\mathbf{u}(x)
\] for almost every $x$. For each $j$, let $X_j$ be a set of full measure on which the integral representation already proved for $\mathbf{u}_j$ holds. Intersecting the sets $X_j$, the set on which the preceding scalar integrals are finite, and the two sets of pointwise convergence gives a single set $X\subseteq M$ of full measure on which all these properties hold simultaneously. If $x\in X$, \[
\begin{split}
&\left|
\mathcal G_s^{\mathbf{E}}\mathbf{u}_j(x)
-
\int_MG_s^{\mathbf{E}}(x,y)\mathbf{u}(y)\,d\lambda_{\mathbf{g}}(y)
\right|_{\mathbf{h}_{\mathbf{E}}}\\
&\qquad\leq
\int_Mg_s(x,y)|\mathbf{u}_j(y)-\mathbf{u}(y)|_{\mathbf{h}_{\mathbf{E}}}\,d\lambda_{\mathbf{g}}(y)
=
\mathcal T_s(|\mathbf{u}_j-\mathbf{u}|_{\mathbf{h}_{\mathbf{E}}})(x),
\end{split}
\] so the left-hand side tends to zero. By uniqueness of the limit, the extension $\mathcal G_s^{\mathbf{E}}\mathbf{u}$ agrees almost everywhere with the kernel integral. Since $\mathcal G_s^{\mathbf{E}}=J_{s,p}^{\mathbf{E}}$ in $L^p(M,\mathbf{E})$, the asserted representation follows. \end{proof}

\subsection{Definition of the fractional scale}

With the potential and its kernel available, we can define real orders of regularity intrinsically. The definition is arranged to agree with integer orders and to be stable under spectral calculus.

\begin{definition}[Bessel Sobolev spaces of positive order] \label{def:sobolev-fraccionario-y-nucleo-calor-espacios-de-sobolev-fraccionarios} Let $(M,\mathbf{g})$ be a connected, complete Riemannian manifold without boundary, and let $\mathbf{E}\to M$ be a smooth vector bundle of finite rank with a bundle metric and a compatible connection. \index{fractional Sobolev space!on a vector bundle} Let $s>0$ and $1\leq p<\infty$. Define \[
H_{L_{\mathbf{E}}}^{s,p}(M,\mathbf{E}):=J_{s,p}^{\mathbf{E}}(L^p(M,\mathbf{E})).
\] If $\mathbf{u}=J_{s,p}^{\mathbf{E}}\mathbf{f}$, set \[
\|\mathbf{u}\|_{H_{L_{\mathbf{E}}}^{s,p}(M,\mathbf{E})}:=\|\mathbf{f}\|_{L^p(M,\mathbf{E})}.
\] Injectivity of $J_{s,p}^{\mathbf{E}}$ ensures that this norm is well defined. For $s=0$, set $H_{L_{\mathbf{E}}}^{0,p}(M,\mathbf{E}):=L^p(M,\mathbf{E})$. \end{definition}

\begin{definition}[Spaces of negative order] \label{def:espacios-bessel-orden-negativo-haces} Let $(M,\mathbf{g})$ be a connected, complete Riemannian manifold without boundary, and let $\mathbf{E}\to M$ be a smooth vector bundle of finite rank with a bundle metric and a compatible connection. Let $s<0$ and $1\leq p<\infty$. Define $H_{L_{\mathbf{E}}}^{s,p}(M,\mathbf{E})$ as the completion of $L^p(M,\mathbf{E})$ with respect to the norm \[
\|\mathbf{u}\|_{H_{L_{\mathbf{E}}}^{s,p}(M,\mathbf{E})}
:=
\|J_{-s,p}^{\mathbf{E}}\mathbf{u}\|_{L^p(M,\mathbf{E})}.
\] \end{definition}

\begin{notation}[Scalar Bessel scale] \label{not:escala-escalar-bessel} On the trivial bundle of rank one, equipped with the trivial connection, the connection Laplacian is the positive realization \(L=\Delta_B=-\Delta\) introduced in the section on the scalar semigroup. For every \(s\in\mathbb R\) and \(1\leq p<\infty\), we write \[
H_L^{s,p}(M):=H_{L_{\mathbf{E}}}^{s,p}(M,\mathbf{E})
\] in this case. Thus \(H_L^{s,p}(M)\) always denotes the intrinsic scale constructed through the heat semigroup and the potentials \((I+L)^{-\displaystyle\frac{s}{2}}\), while the notation without a subscript is reserved for the spaces defined by localization when these are introduced. \end{notation}

\begin{proposition}[Conjugation of the Bessel scale] \label{prop:conjugacion-escala-bessel-haces} Let $(M,\mathbf{g})$ be a connected, complete Riemannian manifold without boundary, and let $\mathbf{E}\to M$ be a smooth complex vector bundle of finite rank with a Hermitian metric and a compatible connection. The conjugate connection $\nabla^{\overline{\mathbf{E}}}$ is compatible with $\mathbf{h}_{\overline{\mathbf{E}}}$, and its Laplacian satisfies \[
 L_{\overline{\mathbf{E}},0}(\overline{\mathbf{u}})=\overline{L_{\mathbf{E},0}\mathbf{u}},
 \qquad \mathbf{u}\in\Gamma_c(\mathbf{E}).
\] For every $1\leq p<\infty$, $t\geq0$, and $r>0$, \begin{align*}
 P_t^{\overline{\mathbf{E}},p}(\overline{\mathbf{u}})
 &=\overline{P_t^{\mathbf{E},p}\mathbf{u}},\\
 J_{r,p}^{\overline{\mathbf{E}}}(\overline{\mathbf{u}})
 &=\overline{J_{r,p}^{\mathbf{E}}\mathbf{u}}.
\end{align*} Consequently, for each $s\in\mathbb R$, the map \begin{equation}
\label{eq:conjugacion-escala-bessel-haces}
 \mathfrak C_{s,p,\mathbf{E}}\colon
 \overline{H_{L_{\mathbf{E}}}^{s,p}(M,\mathbf{E})}
 \longrightarrow
 H_{L_{\overline{\mathbf{E}}}}^{s,p}(M,\overline{\mathbf{E}}),
 \qquad
 \mathfrak C_{s,p,\mathbf{E}}(\overline{\mathbf{u}})=\overline{\mathbf{u}},
\end{equation} is a complex-linear isometric isomorphism. \end{proposition}

\begin{proof} Compatibility of $\nabla^{\overline{\mathbf{E}}}$ and parallelism of $\mathcal R_{\mathbf{E}}$ were proved in Proposition~\ref{prop:riesz-paralelo-conexion-conjugada}. Conjugating the local formula $L_{\mathbf{E},0}\mathbf{u}=-\operatorname{tr}_{\mathbf{g}}(\nabla^2\mathbf{u})$ gives the first identity. Conjugation of sections is an antilinear isometry between the $L^2(M,\mathbf E)$ spaces and between the domains of the energy forms; uniqueness of the associated self-adjoint realization gives $L_{\overline{\mathbf{E}}}\mathfrak C=\mathfrak C L_{\mathbf{E}}$. Spectral calculus, followed by uniqueness of the consistent extensions to $L^p(M,\mathbf E)$, gives the semigroup identity. Integrating this identity in the Laplace formula for $J_{r,p}$ gives the potential identity.

For $s>0$, the norm of the image of $J_{s,p}^{\mathbf{E}}\mathbf{u}$ is the norm of $\overline{\mathbf{u}}$ in $L^p(M,\overline{\mathbf{E}})$, equal to that of $\mathbf{u}$. For $s<0$, the same conclusion follows from the norm defining the completion. Surjectivity follows by conjugating again. Finally, the scalar structure of the conjugate space shows that \eqref{eq:conjugacion-escala-bessel-haces} is complex-linear. \end{proof}

\begin{theorem}[Properties of the Bessel scale] \label{teo:propiedades-escala-bessel-haces} Let $(M,\mathbf{g})$ be a connected, complete Riemannian manifold without boundary, and let $\mathbf{E}\to M$ be a smooth vector bundle of finite rank with a bundle metric and a compatible connection. Let $1\leq p<\infty$. \begin{enumerate}[label=(\alph*)] \item $H_{L_{\mathbf{E}}}^{s,p}(M,\mathbf{E})$ is a Banach space for every $s\in\mathbb R$; \item if $s_1>s_0$, the inclusion \[
H_{L_{\mathbf{E}}}^{s_1,p}(M,\mathbf{E})\hookrightarrow H_{L_{\mathbf{E}}}^{s_0,p}(M,\mathbf{E})
\] is continuous and has dense image; \item for $r>0$ and $s\in\mathbb R$, $J_{r,p}^{\mathbf{E}}$ extends to an isometric isomorphism \[
J_{r,p}^{\mathbf{E}}\colon H_{L_{\mathbf{E}}}^{s,p}(M,\mathbf{E})\longrightarrow H_{L_{\mathbf{E}}}^{s+r,p}(M,\mathbf{E});
\] \item if $s>0$, $J_{s,p}^{\mathbf{E}}(\Gamma_c(\mathbf{E}))$ is dense in $H_{L_{\mathbf{E}}}^{s,p}(M,\mathbf{E})$; \item if $1<p<\infty$, $p'$ is the conjugate exponent, and $s\in\mathbb R$, integral duality, initially defined on $L^p(M,\mathbf{E})\times
L^{p'}(M,\overline{\mathbf{E}})$, extends uniquely to a continuous complex-bilinear pairing \[
 \mathfrak b_{s,p;\mathbf{E}}\colon
 H_{L_{\mathbf{E}}}^{s,p}(M,\mathbf{E})\times
 H_{L_{\overline{\mathbf{E}}}}^{-s,p'}(M,\overline{\mathbf{E}})
 \longrightarrow\mathbb C
\] and induces the complex-linear isometric isomorphism \begin{equation}
\label{eq:dualidad-compleja-escala-bessel-haces}
 H_{L_{\overline{\mathbf{E}}}}^{-s,p'}(M,\overline{\mathbf{E}})
 \longrightarrow
 \bigl(H_{L_{\mathbf{E}}}^{s,p}(M,\mathbf{E})\bigr)',
 \qquad
 \mathbf{w}\longmapsto
 \left[\mathbf{u}\longmapsto
 \mathfrak b_{s,p;\mathbf{E}}(\mathbf{u},\mathbf{w})\right].
\end{equation} Equivalently, using \eqref{eq:conjugacion-escala-bessel-haces}, we obtain \[
 \overline{H_{L_{\mathbf{E}}}^{-s,p'}(M,\mathbf{E})}
 \longrightarrow
 \bigl(H_{L_{\mathbf{E}}}^{s,p}(M,\mathbf{E})\bigr)',
 \qquad
 \overline{\mathbf{v}}\longmapsto
 \left[\mathbf{u}\longmapsto
 \int_M\mathbf{h}_{\mathbf{E}}(\mathbf{u},\mathbf{v})\,d\lambda_{\mathbf{g}}\right].
\] In the last formula, the integral has its literal meaning when $\mathbf{u}\in L^p(M,\mathbf{E})$ and $\mathbf{v}\in L^{p'}(M,\mathbf{E})$; for elements of general orders, it denotes, by definition, the unique continuous extension $\mathfrak b_{s,p;\mathbf{E}}(\mathbf{u},\mathfrak C_{-s,p',\mathbf{E}}(\overline{\mathbf{v}}))$. This convention does not assign pointwise values to a distribution of negative order. \end{enumerate} \end{theorem}

\begin{proof} Parts (a)--(c) are Theorem~\ref{teo:escala-bessel-semigrupo-contraccion} applied to the semigroup $(P_t^{\mathbf{E},p})_{t\geq0}$. For (d), if $\mathbf{u}=J_{s,p}^{\mathbf{E}}\mathbf{f}$, choose $\mathbf{f}_j\in\Gamma_c(\mathbf{E})$ with $\mathbf{f}_j\to \mathbf{f}$ in $L^p(M,\mathbf{E})$. Then \[
\|J_{s,p}^{\mathbf{E}}\mathbf{f}_j-\mathbf{u}\|_{H_{L_{\mathbf{E}}}^{s,p}(M,\mathbf{E})}
=
\|\mathbf{f}_j-\mathbf{f}\|_{L^p(M,\mathbf{E})}
\longrightarrow0.
\]

For duality, self-adjointness on $L^2(M,\mathbf E)$ and consistency of the extensions imply \[
 \int_M\mathbf{h}_{\mathbf{E}}(P_t^{\mathbf{E},p}\mathbf{u},\mathbf{v})\,d\lambda_{\mathbf{g}}
 =\int_M\mathbf{h}_{\mathbf{E}}(\mathbf{u},P_t^{\mathbf{E},p'}\mathbf{v})\,d\lambda_{\mathbf{g}}
\] for $\mathbf{u}$ and $\mathbf{v}$ in dense intersections. The integral defining $J_{r,p}^{\mathbf{E}}$ is dominated by \[
 \frac{t^{\frac{r}{2}-1}e^{-t}}{\Gamma(\frac{r}{2})}
 \|\mathbf{u}\|_{L^p(M,\mathbf{E})}\|\mathbf{v}\|_{L^{p'}(M,\mathbf{E})},
\] which is integrable on $(0,\infty)$. Fubini's theorem for Bochner integrals~\ref{teo:b5-fubini-bochner} gives \begin{equation}
\label{eq:simetria-potencial-bessel-hermitiana}
 \int_M\mathbf{h}_{\mathbf{E}}(J_{r,p}^{\mathbf{E}}\mathbf{u},\mathbf{v})\,d\lambda_{\mathbf{g}}
 =\int_M\mathbf{h}_{\mathbf{E}}(\mathbf{u},J_{r,p'}^{\mathbf{E}}\mathbf{v})\,d\lambda_{\mathbf{g}}.
\end{equation}

Proposition~\ref{prop:dualidad-escala-bessel-consistente}, applied to the bilinear pairing $L^p(M,\mathbf{E})\times L^{p'}(M,\overline{\mathbf{E}})$ induced by $\mathcal R_{\mathbf{E}}$, uses \eqref{eq:simetria-potencial-bessel-hermitiana} and proves that \eqref{eq:dualidad-compleja-escala-bessel-haces} is isometric and surjective. It is complex-linear in $\mathbf{w}$ because $\mathcal R_{\mathbf{E}}\colon\overline{\mathbf{E}}\to \mathbf{E}^*$ is complex-linear. Composition with the linear isomorphism $\mathfrak C_{-s,p',\mathbf{E}}$ gives the equivalent formulation from the conjugate space. In the real case, $\overline{\mathbf{E}}=\mathbf{E}$ and we recover the usual real duality. \end{proof}

\subsection{Spectral characterization in \texorpdfstring{$L^2(M,\mathbf E)$}{L2}}

\begin{theorem}[Spectral characterization] \label{teo:caracterizacion-espectral-bessel-haces} Let $(M,\mathbf{g})$ be a connected, complete Riemannian manifold without boundary, and let $\mathbf{E}\to M$ be a smooth vector bundle of finite rank with a bundle metric and a compatible connection. Let $s>0$. Then \[
J_{s,2}^{\mathbf{E}}=(I+L_{\mathbf{E}})^{-\frac{s}{2}}
\] and \[
H_{L_{\mathbf{E}}}^{s,2}(M,\mathbf{E})
=
\mathcal D\bigl((I+L_{\mathbf{E}})^{\frac{s}{2}}\bigr).
\] For $\mathbf{u}$ in this space, \[
\|\mathbf{u}\|_{H_{L_{\mathbf{E}}}^{s,2}(M,\mathbf{E})}
=
\|(I+L_{\mathbf{E}})^{\frac{s}{2}}\mathbf{u}\|_{L^2(M,\mathbf{E})}.
\] In particular, \[
H_{L_{\mathbf{E}}}^{1,2}(M,\mathbf{E})=W^{1,2}(M,\mathbf{E})
\] and \[
\|\mathbf{u}\|_{H_{L_{\mathbf{E}}}^{1,2}(M,\mathbf{E})}^2
=
\|\mathbf{u}\|_{L^2(M,\mathbf{E})}^2
+
\|\nabla^{\mathbf{E}}\mathbf{u}\|_{L^2(M,T^*M\otimes \mathbf{E})}^2.
\] \end{theorem}

\begin{proof} For $\lambda\geq0$, \[
(1+\lambda)^{-\frac{s}{2}}
=
\frac{1}{\Gamma(\frac{s}{2})}
\int_0^\infty t^{\frac{s}{2}-1}e^{-t}e^{-t\lambda}\,dt.
\] To justify the passage to spectral calculus precisely, fix $\mathbf{u}\in L^2(M,\mathbf{E})$ and, for $0<\varepsilon<R<\infty$, set \[
m_{\varepsilon,R}(\lambda)
:=
\frac{1}{\Gamma(\frac{s}{2})}
\int_\varepsilon^R
t^{\frac{s}{2}-1}e^{-t}e^{-t\lambda}\,dt.
\] On the compact interval $[\varepsilon,R]$, the map $t\mapsto e^{-tL_{\mathbf{E}}}$ is continuous in operator norm: indeed, the spectral calculus in Theorem~\ref{teo:calculo-espectral-autoadjunto-no-negativo} and the scalar fundamental theorem of calculus give \[
\|e^{-tL_{\mathbf{E}}}-e^{-\tau L_{\mathbf{E}}}\|_{\mathcal L(L^2(M,\mathbf{E}))}
\leq
|t-\tau|\sup_{\lambda\in[0,\infty)}\lambda e^{-\varepsilon\lambda}
\leq
\frac{|t-\tau|}{e\varepsilon}.
\] Approximating the integral by Riemann sums and using the rules of the same functional calculus yields \[
\frac{1}{\Gamma(\frac{s}{2})}
\int_\varepsilon^R
t^{\frac{s}{2}-1}e^{-t}e^{-tL_{\mathbf{E}}}\mathbf{u}\,dt
=
m_{\varepsilon,R}(L_{\mathbf{E}})\mathbf{u}.
\] As $\varepsilon\to0^+$ and $R\to\infty$, the left-hand side converges to $J_{s,2}^{\mathbf{E}}\mathbf{u}$ by the Bochner dominated convergence theorem~\ref{teo:b5-convergencia-dominada-bochner}. On the other hand, $m_{\varepsilon,R}(\lambda)\to(1+\lambda)^{-\displaystyle\frac{s}{2}}$ and $0\leq m_{\varepsilon,R}(\lambda)\leq1$. The norm formula in Theorem~\ref{teo:calculo-espectral-autoadjunto-no-negativo}, followed by Lebesgue's dominated convergence theorem~\ref{convergencia dominada}, implies \[
m_{\varepsilon,R}(L_{\mathbf{E}})\mathbf{u}
\longrightarrow
(I+L_{\mathbf{E}})^{-\frac{s}{2}}\mathbf{u}
\quad\text{in }L^2(M,\mathbf{E}).
\] Consequently, $J_{s,2}^{\mathbf{E}}\mathbf{u}=(I+L_{\mathbf{E}})^{-\displaystyle\frac{s}{2}}\mathbf{u}$. Corollary~\ref{cor:potencias-espectrales-operador-no-negativo} identifies the range of $(I+L_{\mathbf{E}})^{-\displaystyle\frac{s}{2}}$ with the domain of $(I+L_{\mathbf{E}})^{\displaystyle\frac{s}{2}}$ and shows that these operators are inverses on the respective spaces.

For $s=1$, Corollary~\ref{cor:dominio-raiz-forma-cerrada} gives \[
\mathcal D\bigl((I+L_{\mathbf{E}})^{\frac{1}{2}}\bigr)=W^{1,2}(M,\mathbf{E})
\] and \[
\|(I+L_{\mathbf{E}})^{\frac{1}{2}}\mathbf{u}\|_{L^2(M,\mathbf{E})}^2
=
\|\mathbf{u}\|_{L^2(M,\mathbf{E})}^2+\mathfrak q_{\mathbf{E}}(\mathbf{u},\mathbf{u}).
\] The definition of $\mathfrak q_{\mathbf{E}}$ completes the proof. \end{proof}

\begin{proposition}[Spectral, dual, and interpolation compatibility in \texorpdfstring{$L^2(M,\mathbf E)$}{L2}] \label{prop:compatibilidad-espectral-dual-interpolacion-bessel-haces} Let $(M,\mathbf{g})$ be a connected, complete Riemannian manifold without boundary, and let $\mathbf{E}\to M$ be a smooth vector bundle of finite rank with a bundle metric $\mathbf h_{\mathbf E}$, Hermitian in the complex case, and a compatible connection. First suppose that $\mathbf{E}$ is complex, and set $B_{\mathbf{E}}:=I+L_{\mathbf{E}}$. Let $s_0<s_1$, $0<\theta<1$, and $s_\theta:=(1-\theta)s_0+\theta s_1$. Then \begin{equation}
\label{eq:interpolacion-hilbert-bessel-haces}
 [H_{L_{\mathbf{E}}}^{s_0,2}(M,\mathbf{E}),H_{L_{\mathbf{E}}}^{s_1,2}(M,\mathbf{E})]_\theta
 =H_{L_{\mathbf{E}}}^{s_\theta,2}(M,\mathbf{E})
\end{equation} as topological spaces, with the norms determined by spectral calculus. This identification is the identity on the intersection of the two endpoints and commutes with Bessel potentials and the conjugation \eqref{eq:conjugacion-escala-bessel-haces}.

Moreover, if $s>0$, $\mathbf{w}\in H_{L_{\overline{\mathbf{E}}}}^{-s,2}(M,\overline{\mathbf{E}})$, \[
 z:=J_{s,2}^{\overline{\mathbf{E}}}\mathbf{w}\in L^2(M,\overline{\mathbf{E}})
\], and $\mathbf{u}\in H_{L_{\mathbf{E}}}^{s,2}(M,\mathbf{E})$, denote by \(\langle \mathbf{u},\mathbf{w}\rangle_{s,2;\mathbf{E}}\) the pairing in Theorem~\ref{teo:propiedades-escala-bessel-haces}. It is given by the spectral formula \begin{equation}
\label{eq:formula-espectral-dualidad-bessel-haces}
 \langle \mathbf{u},\mathbf{w}\rangle_{s,2;\mathbf{E}}
 =
 \int_M
 \mathcal R_{\mathbf{E}}(z)
 \bigl(B_{\mathbf{E}}^{\frac{s}{2}}\mathbf{u}\bigr)
 \,d\lambda_{\mathbf{g}}.
\end{equation} In particular, \[
 H_{L_{\overline{\mathbf{E}}}}^{-s,2}(M,\overline{\mathbf{E}})
 \xrightarrow{\ \cong\ }
 \bigl(\mathcal D(B_{\mathbf{E}}^{\frac{s}{2}}),
       \|B_{\mathbf{E}}^{\frac{s}{2}}\,\cdot\|_{L^2(M,\mathbf{E})}\bigr)'
\] is precisely the identification in \eqref{eq:dualidad-compleja-escala-bessel-haces}; passage to the spectral characterization introduces no second duality convention. When interpolating two orders, the pairing at order $s_\theta$ is the unique continuous extension of the same pairings on the common spectral core. Thus duality, conjugation, and interpolation give compatible identifications.

If $\mathbf{E}$ is real, the same assertions are obtained in the complexification and restricted to the invariant real subspace; in the dual formula, $\overline{\mathbf{E}}$ and $\mathcal R_{\mathbf{E}}$ are omitted and the real bundle inner product is used. \end{proposition}

\begin{proof} First prove \eqref{eq:interpolacion-hilbert-bessel-haces} when $0\leq s_0<s_1$. Write \[
 a_j:=\frac{s_j}{2},\qquad
 a(z):=(1-z)a_0+za_1,\qquad
 a:=a(\theta),\qquad
 \delta:=a_1-a_0.
\] By the preceding theorem, $H_{L_{\mathbf{E}}}^{s_j,2}(M,\mathbf{E})=\mathcal D(B_{\mathbf{E}}^{a_j})$, with norm $\|B_{\mathbf{E}}^{a_j}\,\cdot\|_{L^2(M,\mathbf{E})}$. Spectral calculus gives \begin{equation}
\label{eq:unitariedad-potencias-imaginarias-bessel-haces}
 \|B_{\mathbf{E}}^{it}\mathbf{v}\|_{L^2(M,\mathbf{E})}=\|\mathbf{v}\|_{L^2(M,\mathbf{E})},
 \qquad t\in\mathbb R.
\end{equation}

Let $\mathbf{u}\in\mathcal D(B_{\mathbf{E}}^a)$. For $\varepsilon>0$, the family \[
 F_\varepsilon(z)
 :=e^{\varepsilon(z-\theta)^2}B_{\mathbf{E}}^{a-a(z)}\mathbf{u}
\] is admissible for the complex method. To verify this, set
\[
 w:=B_{\mathbf E}^a\mathbf u,\qquad
 E_N:=\mathbf 1_{[1,N]}(B_{\mathbf E}),\qquad
 F_{\varepsilon,N}(z):=e^{\varepsilon(z-\theta)^2}
                         B_{\mathbf E}^{-a(z)}E_Nw.
\]
Each truncated family is holomorphic with values in all spectral
domains. If
\(C_\varepsilon:=e^{\varepsilon\max\{\theta^2,(1-\theta)^2\}}\), then,
for \(0\leq\operatorname{Re}z\leq1\), spectral calculus gives
\[
 \|F_{\varepsilon,N}(z)-F_\varepsilon(z)\|_{H_{L_{\mathbf E}}^{s_0,2}}
 \leq C_\varepsilon\|(I-E_N)w\|_{L^2},
\]
since \(a(z)=a_0+\delta z\) and \(\|B_{\mathbf E}^{-\delta z}\|\leq1\).
On both boundaries we also obtain
\[
 \sup_{t\in\mathbb R}
 \|F_{\varepsilon,N}(j+it)-F_\varepsilon(j+it)\|_{H_{L_{\mathbf E}}^{s_j,2}}
 \leq C_\varepsilon\|(I-E_N)w\|_{L^2},\qquad j\in\{0,1\}.
\]
The right-hand side tends to zero because \(E_Nw\to w\) in \(L^2\).
These bounds prove continuity on the closed strip, convergence in both
boundary norms, and boundedness in the sum space.
They also justify holomorphy and parameter derivatives. Indeed,
if \(\operatorname{Re}z\geq\eta>0\) and \(k\geq1\), then
\[
 \begin{aligned}
 &\|\partial_z^k\bigl(B_{\mathbf E}^{-a(z)}(I-E_N)w\bigr)
                                      \|_{H_{L_{\mathbf E}}^{s_0,2}}\\
 &\quad\leq\delta^k
       \sup_{\lambda\geq1}(\log\lambda)^k\lambda^{-\delta\eta}
                           \|(I-E_N)w\|_{L^2}
 =\left(\frac{k}{e\eta}\right)^k\|(I-E_N)w\|_{L^2}.
 \end{aligned}
\]
The equality follows by maximizing \(x^ke^{-\delta\eta x}\) for
\(x=\log\lambda\geq0\). Derivatives of the Gaussian factor are bounded
on each compact subset of the strip; Leibniz's rule gives convergence
of all derivatives there. The limit is therefore holomorphic in the
interior. The factor \(e^{-\varepsilon t^2}\) makes its boundary values
tend to zero as \(|t|\to\infty\).

For $j\in\{0,1\}$ and $t\in\mathbb R$, \[
 B_{\mathbf{E}}^{a_j}F_\varepsilon(j+it)
 =e^{\varepsilon(j+it-\theta)^2}B_{\mathbf{E}}^{a-it\delta}\mathbf{u}.
\] Identity \eqref{eq:unitariedad-potencias-imaginarias-bessel-haces} gives \[
 \|F_\varepsilon(j+it)\|_{H_{L_{\mathbf{E}}}^{s_j,2}(M,\mathbf{E})}
 \leq
 e^{\varepsilon(j-\theta)^2}e^{-\varepsilon t^2}
 \|B_{\mathbf{E}}^a\mathbf{u}\|_{L^2(M,\mathbf{E})}.
\] Since $F_\varepsilon(\theta)=\mathbf{u}$, we let $\varepsilon\to0^+$ and obtain the embedding of $H_{L_{\mathbf{E}}}^{s_\theta,2}(M,\mathbf{E})$ into the left-hand side of \eqref{eq:interpolacion-hilbert-bessel-haces}.

For the reverse embedding, let $F$ be an admissible family with $\mathbf{u}=F(\theta)$, and let \[
 \mathscr C_{\mathbf{E}}
 :=\bigcup_{N\in\mathbb N}
 \mathbf 1_{[1,N]}(B_{\mathbf{E}})L^2(M,\mathbf{E}).
\] This subspace is a core for every nonnegative real power of $B_{\mathbf{E}}$. For $\mathbf{v}\in\mathscr C_{\mathbf{E}}$, the function \[
 \Phi_v(z)
 :=
 \int_M \mathbf{h}_{\mathbf{E}}\bigl(F(z),B_{\mathbf{E}}^{a(\overline z)}\mathbf{v}\bigr)
 \,d\lambda_{\mathbf{g}}
\] is holomorphic: the second entry is antiholomorphic, and $\mathbf{h}_{\mathbf{E}}$ is antilinear in that entry. On the boundary $z=j+it$, self-adjointness and \eqref{eq:unitariedad-potencias-imaginarias-bessel-haces} imply \[
 |\Phi_v(j+it)|
 \leq
 \|F(j+it)\|_{H_{L_{\mathbf{E}}}^{s_j,2}(M,\mathbf{E})}
 \|B_{\mathbf{E}}^{-a_j}B_{\mathbf{E}}^{a_j-it\delta}\mathbf{v}\|_{L^2(M,\mathbf{E})}
 =
 \|F(j+it)\|_{H_{L_{\mathbf{E}}}^{s_j,2}(M,\mathbf{E})}\|\mathbf{v}\|_{L^2(M,\mathbf{E})}.
\] The three-lines lemma~\ref{lem:tres-lineas-hadamard} gives \[
 \left|
 \int_M\mathbf{h}_{\mathbf{E}}(\mathbf{u},B_{\mathbf{E}}^a\mathbf{v})\,d\lambda_{\mathbf{g}}
 \right|
 \leq
 \|F\|_{\mathscr F}\|\mathbf{v}\|_{L^2(M,\mathbf{E})}.
\] Since $\mathscr C_{\mathbf{E}}$ is a core for $B_{\mathbf{E}}^a$, conjugating the preceding inequality gives a bounded linear functional \[
 \mathbf{v}\longmapsto
 \int_M\mathbf{h}_{\mathbf{E}}(B_{\mathbf{E}}^a\mathbf{v},\mathbf{u})\,d\lambda_{\mathbf{g}}
\] on $\mathscr C_{\mathbf{E}}$ with respect to the norm of $L^2(M,\mathbf{E})$. This extends to all of $L^2(M,\mathbf{E})$, and the Riesz representation theorem~\ref{teo:representacion-riesz-hilbert} provides $\mathbf{g}\in L^2(M,\mathbf{E})$ such that \[
 \int_M\mathbf{h}_{\mathbf{E}}(B_{\mathbf{E}}^a\mathbf{v},\mathbf{u})\,d\lambda_{\mathbf{g}}
 =\int_M\mathbf{h}_{\mathbf{E}}(\mathbf{v},\mathbf{g})\,d\lambda_{\mathbf{g}},
 \qquad \mathbf{v}\in\mathscr C_{\mathbf{E}}.
\] Since $\mathscr C_{\mathbf{E}}$ is a core for $B_{\mathbf{E}}^a$, the identity extends to all of $v\in\mathcal D(B_{\mathbf{E}}^a)$. By the definition of the adjoint and self-adjointness of $B_{\mathbf{E}}^a$, this proves that $\mathbf{u}\in\mathcal D(B_{\mathbf{E}}^a)$, $B_{\mathbf{E}}^a\mathbf{u}=\mathbf{g}$, and \(
\|B_{\mathbf{E}}^a\mathbf{u}\|_{L^2(M,\mathbf{E})}\leq\|F\|_{\mathscr F}
\). Taking the infimum over $F$ gives the reverse embedding.

For arbitrary orders $s_0<s_1$, choose $\displaystyle r>\displaystyle\max\{0,-s_0\}$. By Theorem~\ref{teo:propiedades-escala-bessel-haces}, the same potential $J_{r,2}^{\mathbf{E}}$ is a simultaneous isometric isomorphism \[
 H_{L_{\mathbf{E}}}^{s_j,2}(M,\mathbf{E})\longrightarrow H_{L_{\mathbf{E}}}^{s_j+r,2}(M,\mathbf{E}),
 \qquad j\in\{0,1\},
\] and also at the intermediate order. The complex interpolation theorem for operators~\ref{teo:interpolacion-compleja-operadores}, applied to $J_{r,2}^{\mathbf{E}}$ and its inverse on these compatible pairs, reduces the general case to the one already proved.

Now prove \eqref{eq:formula-espectral-dualidad-bessel-haces}. Set $\mathbf{f}:=B_{\mathbf{E}}^{\displaystyle\frac{s}{2}}\mathbf{u}$, so that $\mathbf{u}=J_{s,2}^{\mathbf{E}}\mathbf{f}$. The bilinear symmetry of the potentials, expressed between $\mathbf{E}$ and $\overline{\mathbf{E}}$, gives \[
 \langle J_{s,2}^{\mathbf{E}}\mathbf{f},\mathbf{w}\rangle_{s,2;\mathbf{E}}
 =
 \int_M\mathcal R_{\mathbf{E}}(J_{s,2}^{\overline{\mathbf{E}}}\mathbf{w})(\mathbf{f})
 \,d\lambda_{\mathbf{g}},
\] which is exactly the asserted formula. The operators $J_{r,2}$, spectral powers of real exponent, and conjugation commute because their scalar multipliers are real. Finally, spectral truncations form a dense subspace at every order; on this subspace, all the preceding identifications are the same algebraic formulas. Uniqueness of continuous extensions proves the final compatibility assertion. \end{proof}

\begin{remark}[Complex interpolation outside $L^2(M,\mathbf E)$] \label{obs:alcance-interpolacion-escala-bessel-semigrupo} The Hilbert space proof uses the fact that $(I+L_{\mathbf E})^{it}$ is unitary on $L^2(M,\mathbf E)$ for $t\in\mathbb R$. For $1<p<\infty$, the argument requires a bound for these powers on $L^p(M,\mathbf E)$; heat contractivity alone does not provide this bound for an abstract semigroup. Lemma~\ref{lem:potencias-imaginarias-bochner-haz-completo} will prove it for the covariant Laplacian on a complete manifold, and Theorem~\ref{teo:interpolacion-compleja-bessel-laplaciano-completo} will establish the corresponding interpolation result. Comparison with localized norms additionally requires bounded geometry of the base and bundle; it is established in Theorem~\ref{teo:bessel-calor-triebel-local-haces-geometria-acotada}. \end{remark}

\begin{corollary} \label{cor:orden-par-dominio-generador-Lp} Let $(M,\mathbf{g})$ be a connected, complete Riemannian manifold without boundary, and let $\mathbf{E}\to M$ be a smooth vector bundle of finite rank with a bundle metric and a compatible connection. If $m\in\mathbb N$ and $1\leq p<\infty$, then \[
H_{L_{\mathbf{E}}}^{2m,p}(M,\mathbf{E})
=
\mathcal D\bigl((I+L_{\mathbf{E},p})^m\bigr)
\] and \[
\|\mathbf{u}\|_{H_{L_{\mathbf{E}}}^{2m,p}(M,\mathbf{E})}
=
\|(I+L_{\mathbf{E},p})^m\mathbf{u}\|_{L^p(M,\mathbf{E})}.
\] \end{corollary}

\begin{proof} Proposition~\ref{prop:propiedades-potenciales-bessel-haces} gives \[
J_{2m,p}^{\mathbf{E}}=(I+L_{\mathbf{E},p})^{-m}.
\] The range of this iterated resolvent is the domain of $(I+L_{\mathbf{E},p})^m$, and the range norm agrees with the indicated norm. \end{proof}

\begin{corollary}[Comparison on a closed manifold] \label{cor:bessel-sobolev-entero-variedad-cerrada} If $M$ is closed, then, for every $m\in\mathbb N_0$, \[
H_{L_{\mathbf{E}}}^{m,2}(M,\mathbf{E})=H^m(M,\mathbf{E})
\] and there exist constants $c_m,C_m>0$, depending only on $m$ and the fixed data $(M,g,\mathbf{E},\mathbf{h}_{\mathbf{E}},\nabla^{\mathbf{E}})$, such that \[
c_m\|\mathbf{u}\|_{H^m(M,\mathbf{E})}
\leq\|\mathbf{u}\|_{H_{L_{\mathbf{E}}}^{m,2}(M,\mathbf{E})}
\leq C_m\|\mathbf{u}\|_{H^m(M,\mathbf{E})},
\qquad \mathbf{u}\in H^m(M,\mathbf{E}).
\] \end{corollary}

\begin{proof} For $m=0$, the two norms agree. Let $m=2q$ with $q\geq1$. Corollary~\ref{cor:potencias-laplaciano-bochner-dominios-sobolev} gives \[
\mathcal D(L_{\mathbf{E}}^q)=H^{2q}(M,\mathbf{E})
\] and constants $a_q,A_q>0$ for which \[
a_q\|\mathbf{u}\|_{H^{2q}(M,\mathbf{E})}
\leq\|\mathbf{u}\|_{L^2(M,\mathbf{E})}+\|L_{\mathbf{E}}^q\mathbf{u}\|_{L^2(M,\mathbf{E})}
\leq A_q\|\mathbf{u}\|_{H^{2q}(M,\mathbf{E})}.
\] For $\lambda\geq0$, \[
1+\lambda^{2q}
\leq(1+\lambda)^{2q}
\leq2^{2q-1}(1+\lambda^{2q}).
\] Spectral calculus, together with \[
\bigl(\|\mathbf{u}\|_{L^2(M,\mathbf{E})}^2+\|L_{\mathbf{E}}^q\mathbf{u}\|_{L^2(M,\mathbf{E})}^2\bigr)^{\frac{1}{2}}
\leq\|\mathbf{u}\|_{L^2(M,\mathbf{E})}+\|L_{\mathbf{E}}^q\mathbf{u}\|_{L^2(M,\mathbf{E})}
\leq\sqrt2\,
\bigl(\|\mathbf{u}\|_{L^2(M,\mathbf{E})}^2+\|L_{\mathbf{E}}^q\mathbf{u}\|_{L^2(M,\mathbf{E})}^2\bigr)^{\frac{1}{2}},
\] proves the asserted inequality with $c_{2q}=\displaystyle\frac{a_q}{\sqrt2}$ and $C_{2q}=2^{q-\displaystyle\frac{1}{2}}A_q$.

For $m=2q+1$, spectral calculus gives \[
\mathcal D\bigl((I+L_{\mathbf{E}})^{q+\frac{1}{2}}\bigr)
=
\left\{\mathbf{u}\in\mathcal D\bigl((I+L_{\mathbf{E}})^q\bigr)\middle|(I+L_{\mathbf{E}})^q\mathbf{u}\in\mathcal D\bigl((I+L_{\mathbf{E}})^{\frac{1}{2}}\bigr)\right\}.
\] Theorem~\ref{teo:caracterizacion-espectral-bessel-haces} identifies the last domain with $W^{1,2}(M,\mathbf{E})$. Set \[
\mathcal N_q(\mathbf{u})^2
:=
\displaystyle\sum_{j=0}^{q}\|L_{\mathbf{E}}^j\mathbf{u}\|_{L^2(M,\mathbf{E})}^2
+
\|\nabla^{\mathbf{E}}L_{\mathbf{E}}^q\mathbf{u}\|_{L^2(M,T^*M\otimes \mathbf{E})}^2.
\] The identity $\|\nabla^{\mathbf{E}}L_{\mathbf{E}}^j\mathbf{u}\|_{L^2(M,T^*M\otimes \mathbf{E})}^2
=\langle L_{\mathbf{E}}^{j+1}\mathbf{u},L_{\mathbf{E}}^j\mathbf{u}\rangle_{L^2(M,\mathbf{E})}$ and Cauchy--Schwarz give \[
\|\nabla^{\mathbf{E}}L_{\mathbf{E}}^j\mathbf{u}\|_{L^2(M,T^*M\otimes \mathbf{E})}^2
\leq\frac{1}{2}\bigl(\|L_{\mathbf{E}}^{j+1}\mathbf{u}\|_{L^2(M,\mathbf{E})}^2
+\|L_{\mathbf{E}}^j\mathbf{u}\|_{L^2(M,\mathbf{E})}^2\bigr),
\qquad 0\leq j<q.
\] If $\|\mathbf{u}\|_{\Delta_B,2q+1}$ denotes the norm in Theorem~\ref{teo:sobolev_cerrada_riemanniana}, these inequalities imply \[
\mathcal N_q(\mathbf{u})^2
\leq\|\mathbf{u}\|_{\Delta_B,2q+1}^2
\leq(q+1)\mathcal N_q(\mathbf{u})^2.
\] The cited theorem provides $\alpha_q,b_q>0$, depending only on $q$ and the fixed geometric data, such that \[
\alpha_q\|\mathbf{u}\|_{H^{2q+1}(M,\mathbf{E})}^2
\leq\|\mathbf{u}\|_{\Delta_B,2q+1}^2
\leq b_q\|\mathbf{u}\|_{H^{2q+1}(M,\mathbf{E})}^2.
\] On the other hand, for $\lambda\geq0$, \[
\frac{1}{q+2}
\left(\sum_{j=0}^{q}\lambda^{2j}+\lambda^{2q+1}\right)
\leq(1+\lambda)^{2q+1}
\leq2^{2q}
\left(\sum_{j=0}^{q}\lambda^{2j}+\lambda^{2q+1}\right).
\] Spectral calculus turns this polynomial comparison into \[
\frac{1}{\sqrt{q+2}}\mathcal N_q(\mathbf{u})
\leq\|(I+L_{\mathbf{E}})^{q+\frac{1}{2}}\mathbf{u}\|_{L^2(M,\mathbf{E})}
\leq2^q\mathcal N_q(\mathbf{u}).
\] Consequently, we may take \[
c_{2q+1}:=\sqrt{\frac{\alpha_q}{(q+1)(q+2)}},
\qquad
C_{2q+1}:=2^q\sqrt{b_q}.
\] Moreover, $\Gamma(\mathbf{E})$ is dense in both norms: it is dense in $H^{2q+1}(M,\mathbf{E})$ by definition and, since $M$ is compact, it is also dense in $H_{L_{\mathbf{E}}}^{2q+1,2}(M,\mathbf{E})$. Indeed, spectral calculus gives \[
J_{2q+1,2}^{\mathbf{E}}(\Gamma(\mathbf{E}))
\subseteq\bigcap_{k\in\mathbb N}\mathcal D(L_{\mathbf{E}}^k)=\Gamma(\mathbf{E}),
\] while part~(d) of Theorem~\ref{teo:propiedades-escala-bessel-haces} states that the subset on the left-hand side is dense in $H_{L_{\mathbf{E}}}^{2q+1,2}(M,\mathbf{E})$. Thus the inequalities extend to the completions and simultaneously prove equality of the two spaces. \end{proof}

\begin{remark}[Comparison beyond the compact case] \label{obs:comparacion-bessel-sobolev-covariante} The equality $H_{L_{\mathbf{E}}}^{1,2}(M,\mathbf{E})=W^{1,2}(M,\mathbf{E})$ holds on every complete manifold. For higher orders and for $p\neq2$, identifying the Bessel scale with the spaces defined by covariant derivatives requires global elliptic estimates. The chapter on elliptic operators provides these estimates on closed manifolds, while the chapter on bounded geometry explains the uniform hypotheses needed in the noncompact case. Completeness alone guarantees the existence of the Bessel scale, but not global equivalence with a norm formed from finitely many covariant derivatives. \end{remark}

\begin{remark}[Choice of operator] \label{obs:eleccion-operador-fraccionarios-haces} The construction can be carried out with other self-adjoint elliptic operators bounded from below, provided that their semigroups admit consistent extensions to $L^p(M,\mathbf{E})$. In this chapter, we chose the connection Laplacian because it is determined by $(g,\mathbf{h}_{\mathbf{E}},\nabla^{\mathbf{E}})$, its quadratic form is exactly the energy of the connection, and Kato's inequality gives a direct comparison with the scalar semigroup. \end{remark}

\section{Laplacian eigenvalues on compact manifolds}
\label{sec:espectro-laplaciano-compacto}

Rellich compactness turns the self-adjoint realization of the Laplacian
into a discrete spectral problem. We first prove this step for a closed
form, then identify eigenfunctions, the Rayleigh quotient, and the heat
expansion. In this section, the nonnegative scalar operator is
$L=-\Delta_{\mathbf g}$: its eigenfunction equation is
$-\Delta_{\mathbf g}\phi=\lambda\phi$.
Background for this development is available in
\cite{Chavel1984Eigenvalues,Eichhorn2007,Borthwick2020Spectral}.

\begin{definition}[Eigenvalue and resolvent]
\label{def:eigenvalor-resolvente-laplaciano}
Let $A$ be a closed, densely defined operator on a complex Hilbert space
$H$. A number $\lambda\in\mathbb C$ is an \textit{eigenvalue} if there
exists $u\in\mathcal D(A)\setminus\{0\}$ with $Au=\lambda u$. Its
eigenspace is $\ker(A-\lambda I)$, and its geometric multiplicity is the
dimension of that space. Write $\sigma_{\mathrm p}(A)$ for the
set of eigenvalues. The resolvent set $\rho(A)$ consists of the
$z\in\mathbb C$ for which $A-zI$ is bijective from $\mathcal D(A)$ onto
$H$ and its inverse is bounded on $H$; the spectrum is
$\sigma(A):=\mathbb C\setminus\rho(A)$. The inverse is also continuous
into $\mathcal D(A)$ with its graph norm, since
$A(A-zI)^{-1}=I+z(A-zI)^{-1}$.
For $A=L_{\mathbf E}$ we speak of eigensections; in the scalar case,
of eigenfunctions.
\end{definition}

\begin{lemma}[Compact form and eigenvector basis]
\label{lem:criterio-forma-compacta-eigenbase}
Let $\mathfrak q$ be a densely defined, closed, nonnegative
sesquilinear form on an infinite-dimensional Hilbert space $H$.
Let $A$ be its associated self-adjoint operator, and set
\[
 V:=\mathcal D(\mathfrak q),\qquad
 \|u\|_V^2:=\|u\|_H^2+\mathfrak q(u,u).
\]
If the inclusion $V\hookrightarrow H$ is compact, then $(I+A)^{-1}$
is compact. There exist an orthonormal basis $(\phi_j)_{j\geq1}$ of
$H$ and numbers
\[
 0\leq\lambda_1\leq\lambda_2\leq\cdots,\qquad
 \lambda_j\longrightarrow+\infty,\qquad A\phi_j=\lambda_j\phi_j.
\]
The spectrum consists precisely of these eigenvalues, repeated according
to their finite multiplicities. If $c_j=\langle u,\phi_j\rangle_H$, then
\begin{align}
 u\in\mathcal D(A)
 &\ \Longleftrightarrow\
 \sum_{j=1}^{\infty}\lambda_j^2|c_j|^2<\infty,
 \label{eq:dominio-operador-serie-eigenbase}\\
 Au&=\sum_{j=1}^{\infty}\lambda_jc_j\phi_j,\notag\\
 u\in V
 &\ \Longleftrightarrow\
 \sum_{j=1}^{\infty}(1+\lambda_j)|c_j|^2<\infty,
 \label{eq:dominio-forma-serie-eigenbase}\\
 \mathfrak q(u,u)&=\sum_{j=1}^{\infty}\lambda_j|c_j|^2.\notag
\end{align}
The first vector series converges in $H$. For $u\in V$, the partial
sums $\sum_{j=1}^Nc_j\phi_j$ converge to $u$ in $V$.
\end{lemma}

\begin{proof}
For $f\in H$, let $u=(I+A)^{-1}f$. The form representation gives
\[
 \|u\|_V^2=\langle f,u\rangle_H
 \leq\|f\|_H\|u\|_H\leq\|f\|_H\|u\|_V.
\]
Thus the resolvent is continuous from $H$ to $V$, and composing it
with the compact inclusion gives a compact operator on $H$. It is
self-adjoint, injective, and positive. Its range is dense because the
orthogonal complement of its range is its kernel.
Theorem~\ref{teo:espectral-compacto-autoadjunto}, proved in the appendix,
provides an orthonormal basis of resolvent eigenvectors with
eigenvalues $0<\mu_j\leq1$ tending to zero. The identity
$(I+A)^{-1}\phi_j=\mu_j\phi_j$ implies $\phi_j\in\mathcal D(A)$ and
$A\phi_j=(\mu_j^{-1}-1)\phi_j$. Ordering
$\lambda_j:=\mu_j^{-1}-1$ yields the stated sequence.

If $u\in\mathcal D(A)$, self-adjointness gives
$\langle Au,\phi_j\rangle_H=\lambda_jc_j$.
Parseval proves necessity in
\eqref{eq:dominio-operador-serie-eigenbase}. Conversely, if that sum
is finite, the partial sums of $u$ and of its image under $A$
converge in $H$. Since $A$ is closed, $u\in\mathcal D(A)$ and its
image is the indicated series.

The vectors $v_j:=(1+\lambda_j)^{-\displaystyle\frac{1}{2}}\phi_j$ are orthonormal in $V$,
since $\mathfrak q(w,\phi_j)=\langle w,A\phi_j\rangle_H$ for $w\in V$.
If $w$ is orthogonal to every $v_j$ in $V$, then
\[
 0=\langle w,v_j\rangle_V
   =(1+\lambda_j)^{\frac{1}{2}}\langle w,\phi_j\rangle_H
\]
for every $j$, and hence $w=0$. They form an orthonormal basis of $V$.
Parseval in $V$ proves \eqref{eq:dominio-forma-serie-eigenbase} and
convergence of the partial sums in the form norm.

Finally, if $z$ differs from every $\lambda_j$, then
$\inf_j|\lambda_j-z|>0$ and
$\sup_j\displaystyle\frac{\lambda_j}{|\lambda_j-z|}<\infty$, because $\lambda_j\to+\infty$.
The map
\[
 f\longmapsto\sum_{j=1}^{\infty}
          \frac{\langle f,\phi_j\rangle_H}{\lambda_j-z}\phi_j
\]
is bounded on $H$, takes values in $\mathcal D(A)$ by
\eqref{eq:dominio-operador-serie-eigenbase}, and is the inverse of
$A-zI$. This identifies the entire spectrum. Finite multiplicity
follows from that of the resolvent eigenvalue $\mu_j$.
\end{proof}

\begin{theorem}[Laplacian spectrum on a closed manifold]
\label{teo:espectro-laplaciano-cerrada-eigenbase}
Suppose that $M$ is closed, $\dim M\geq1$, and
$\operatorname{rank}\mathbf E\geq1$. The operator $L_{\mathbf E}$
has compact resolvent and an orthonormal basis of smooth eigensections
$(\boldsymbol\phi_j)_{j\geq1}$, with eigenvalues
$0\leq\lambda_1\leq\lambda_2\leq\cdots\to+\infty$.
Its kernel is the space of parallel sections.
Formulas \eqref{eq:dominio-operador-serie-eigenbase} and
\eqref{eq:dominio-forma-serie-eigenbase} describe, respectively,
$H^2(M,\mathbf E)$ and $H^1(M,\mathbf E)$.

In the scalar case,
\[
 \lambda_1=0,\qquad
 \phi_1=\lambda_{\mathbf g}(M)^{-\frac{1}{2}},\qquad
 \lambda_2>0.
\]
The zero eigenvalue is simple because $M$ is connected.
\end{theorem}

\begin{proof}
The form of $L_{\mathbf E}$ has domain $H^1(M,\mathbf E)$ and equals
$\mathfrak q_{\mathbf E}(\mathbf u,\mathbf u)=\|\nabla\mathbf u\|_{L^2}^2$.
The inclusion $H^1(M,\mathbf E)\hookrightarrow L^2(M,\mathbf E)$
is compact by Corollary~\ref{cor:rellich-descenso-un-orden-haces}.
The preceding lemma gives compact resolvent, the basis, and the series
domains; the operator domain is $H^2(M,\mathbf E)$ by
Corollary~\ref{cor:caso-cerrado-laplaciano-conexion-capitulo-fraccionario}.
The equation $L_{\mathbf E}\boldsymbol\phi_j=\lambda_j\boldsymbol\phi_j$
and elliptic regularity, applied successively in charts, give
$\boldsymbol\phi_j\in H^m(M,\mathbf E)$ for every integer $m$.
Sobolev embeddings imply smoothness.

If $\mathbf u\in\ker L_{\mathbf E}$, then
$\|\nabla\mathbf u\|_{L^2}^2=0$; the same elliptic regularity gives
$\nabla\mathbf u=0$ pointwise. A parallel section belongs to the
kernel by the Laplacian formula. For functions, vanishing gradient
is equivalent to being constant on the connected manifold.
Constants are integrable because $M$ is compact. This proves
simplicity of zero and the normalization of $\phi_1$.
Since eigenvalues have no finite accumulation, the next one is
strictly positive.
\end{proof}

\begin{proposition}[Rayleigh quotient and minimax principles]
\label{prop:rayleigh-minimax-forma-compacta}
Under the hypotheses of Lemma~\ref{lem:criterio-forma-compacta-eigenbase},
define $\mathcal R(u):=\displaystyle\frac{\mathfrak q(u,u)}{\|u\|_H^2}$ for
$u\in V\setminus\{0\}$. For every $k\geq1$,
\begin{align}
 \lambda_k
 &=\min_{\substack{W\subset V\\\dim W=k}}
       \ \max_{u\in W\setminus\{0\}}\mathcal R(u),
 \label{eq:minimax-eigenvalores-compactos}\\
 &=\max_{\substack{S\subset H\\\dim S=k-1}}
       \ \inf_{\substack{u\in V\cap S^\perp\\u\neq0}}\mathcal R(u).
 \label{eq:maximin-eigenvalores-compactos}
\end{align}
If $m_0:=\dim\ker A$ and $\lambda_+:=\lambda_{m_0+1}$, then
\[
 \mathfrak q(u,u)\geq\lambda_+\|u\|_H^2,
 \qquad u\in V\cap(\ker A)^\perp.
\]
The constant $\lambda_+$ is optimal. For $A=-\Delta_{\mathbf g}$
on the connected closed manifold, this gives
\begin{equation}
 \|u-\overline u\|_{L^2(M)}^2
 \leq\lambda_2^{-1}\|du\|_{L^2(M,T^*M)}^2,\qquad
 \overline u:=\frac{1}{\lambda_{\mathbf g}(M)}
                   \int_Mu\,d\lambda_{\mathbf g}.
 \label{eq:poincare-eigenvalor-compacto}
\end{equation}
\end{proposition}

\begin{proof}
The form and norm expansions give
\[
 \mathcal R(u)=
 \frac{\sum_{j=1}^{\infty}\lambda_j|c_j|^2}{\sum_{j=1}^{\infty}|c_j|^2}.
\]
If $u$ is orthogonal to $\phi_1,\ldots,\phi_{k-1}$, this expression
is at least $\lambda_k$. Every $k$-dimensional subspace $W\subset V$
contains a nonzero vector satisfying those $k-1$ orthogonality
conditions. Thus the maximum of $\mathcal R$ on $W\setminus\{0\}$
is at least $\lambda_k$. The maximum exists because the unit sphere
of $W$ is compact and $\mathfrak q$ is continuous on that
finite-dimensional space.
For $W=\operatorname{span}\{\phi_1,\ldots,\phi_k\}$, the maximum
equals $\lambda_k$, proving
\eqref{eq:minimax-eigenvalores-compactos}.

If $S\subset H$ has dimension $k-1$, there is a nonzero vector in
$\operatorname{span}\{\phi_1,\ldots,\phi_k\}\cap S^\perp$, and its
quotient is at most $\lambda_k$. Hence the infimum in
\eqref{eq:maximin-eigenvalores-compactos} is at most $\lambda_k$.
Choosing $S=\operatorname{span}\{\phi_1,\ldots,\phi_{k-1}\}$ gives
$\lambda_k$ and proves the second formula.

On the orthogonal complement of the kernel, all coefficients with
$\lambda_j=0$ vanish, and the same sum gives the inequality with
$\lambda_+$. Equality holds for an eigenvector of $\lambda_+$.
In the scalar case, the kernel projection is the constant
function $\overline u$, and $d(u-\overline u)=du$; this yields
\eqref{eq:poincare-eigenvalor-compacto}.
\end{proof}

\begin{theorem}[Lichnerowicz bound]
\label{teo:lichnerowicz-primer-eigenvalor}
Let $(M,\mathbf g)$ be closed and connected, of dimension $n\geq2$.
If $\operatorname{Ric}_{\mathbf g}\geq(n-1)\kappa\mathbf g$
with $\kappa>0$, every positive scalar Laplacian eigenvalue
satisfies $\lambda\geq n\kappa$.
In particular, $\lambda_2\geq n\kappa$.
\end{theorem}

\begin{proof}
Let $\phi$ be a nonzero real eigenfunction of
$L=-\Delta_{\mathbf g}$ with eigenvalue $\lambda>0$.
It can be chosen real because the operator preserves real functions.
Integration by parts gives
$\|d\phi\|_{L^2}^2=\lambda\|\phi\|_{L^2}^2$.
The Hodge Laplacian on one-forms satisfies
$\Delta_1d\phi=dL\phi=\lambda d\phi$:
expand $d\delta+\delta d$ and use $d^2=0$ to obtain commutation.
The Weitzenb\"ock identity from
Theorem~\ref{teo:operadores-diferenciales-en-haces-identidad-de-weitzenbock},
integrated against $d\phi$, gives
\[
 \lambda\|d\phi\|_{L^2}^2
 =\|\nabla d\phi\|_{L^2}^2+
       \int_M\operatorname{Ric}_{\mathbf g}
             (\operatorname{grad}\phi,\operatorname{grad}\phi)\,
             d\lambda_{\mathbf g}.
\]
For a symmetric matrix $H$, Cauchy--Schwarz gives
$(\operatorname{tr}H)^2\leq n|H|^2$.
Apply this to the Hessian of $\phi$, whose trace is
$\Delta_{\mathbf g}\phi=-\lambda\phi$.
Together with the Ricci bound, this yields
\[
 \lambda^2\|\phi\|_{L^2}^2
 \geq\frac{\lambda^2}{n}\|\phi\|_{L^2}^2
      +(n-1)\kappa\lambda\|\phi\|_{L^2}^2.
\]
Division by $(n-1)\lambda\|\phi\|_{L^2}^2>0$ gives
$\lambda\geq n\kappa$.
\end{proof}

\begin{corollary}[Spectral heat series]
\label{cor:calor-series-espectrales-compactas}
Under the hypotheses of
Theorem~\ref{teo:espectro-laplaciano-cerrada-eigenbase}, for
$\mathbf u\in L^2(M,\mathbf E)$,
\begin{align}
 P_t^{\mathbf E}\mathbf u
 &=\sum_{j=1}^{\infty}e^{-t\lambda_j}
          \langle\mathbf u,\boldsymbol\phi_j\rangle_{L^2}
          \boldsymbol\phi_j,
 \label{eq:calor-serie-eigenfunciones}\\
 \mathbf K_t^{\mathbf E}(x,y)
 &=\sum_{j=1}^{\infty}e^{-t\lambda_j}
          \boldsymbol\phi_j(x)\otimes\boldsymbol\phi_j(y)^*.
 \label{eq:nucleo-calor-serie-eigenfunciones}
\end{align}
The first series converges in $L^2$ for $t\geq0$.
For every $\tau>0$, both series and all their spatial and temporal
derivatives converge uniformly for $t\geq\tau$.
In \eqref{eq:nucleo-calor-serie-eigenfunciones}, the tensor term sends
$v\in\mathbf E_y$ to
$\langle v,\boldsymbol\phi_j(y)\rangle_{\mathbf h_{\mathbf E}}
\boldsymbol\phi_j(x)$.

If $\Pi_0$ is the orthogonal projection onto $\ker L_{\mathbf E}$, then
\[
 \|P_t^{\mathbf E}\mathbf u-\Pi_0\mathbf u\|_{L^2}
 \leq e^{-t\lambda_+}\|\mathbf u-\Pi_0\mathbf u\|_{L^2}.
\]
In particular, scalar heat converges to the mean of its initial datum.
\end{corollary}

\begin{proof}
Functional calculus on each eigenvector gives
$P_t^{\mathbf E}\boldsymbol\phi_j=e^{-t\lambda_j}\boldsymbol\phi_j$.
Since the semigroup is bounded and the basis is complete, this proves
\eqref{eq:calor-serie-eigenfunciones} in $L^2$.
Let $Q_N$ be the projections onto the first $N$ eigenvectors.
For integers $a,b\geq0$ and $t\geq\tau>0$,
\[
 \|(I+L_{\mathbf E})^a L_{\mathbf E}^b
       P_t^{\mathbf E}(I-Q_N)\mathbf u\|_{L^2}
 \leq C_{a,b,\tau}\|(I-Q_N)\mathbf u\|_{L^2},
\]
where
$C_{a,b,\tau}:=\sup_{\lambda\geq0}(1+\lambda)^a\lambda^b
e^{-\tau\lambda}<\infty$.
The equivalence between the domain of $(I+L_{\mathbf E})^a$ and
$H^{2a}$, followed by $H^{2a}\hookrightarrow C^r$ for
$2a>r+\displaystyle\frac{\dim M}{2}$, proves uniform convergence of the derivatives
of the first series.
A temporal derivative of order $b$ introduces $(-\lambda_j)^b$;
the same estimate justifies termwise differentiation.

On $\operatorname{Hom}(\mathbf E_y,\mathbf E_x)$ we use the tensor
metric: for an orthonormal basis $(e_a)$ of $\mathbf E_y$,
$|T|^2=\sum_a|Te_a|^2$, independently of that basis.
The previously constructed smooth kernel $\mathbf K_t^{\mathbf E}$ is square
integrable on $M\times M$. Fubini and Parseval in both variables give
an orthonormal basis of $L^2$ of the bundle
$\mathbf E\boxtimes\mathbf E^*$ consisting of
$\boldsymbol\phi_i(x)\otimes\boldsymbol\phi_j(y)^*$.
To check completeness, if a kernel $S$ is orthogonal to all these
elements, set
\[
 a_j(x):=\int_M S(x,y)\boldsymbol\phi_j(y)\,d\lambda_{\mathbf g}(y).
\]
Fubini and Cauchy--Schwarz give $a_j\in L^2(M,\mathbf E)$.
Its inner product with every $\boldsymbol\phi_i$ is zero, so $a_j=0$.
Outside a single null set this holds for all $j$ simultaneously.
In each fiber $\mathbf E_x$, Parseval in the variable $y$, applied
to the components of $S(x,\cdot)$, gives
$S(x,\cdot)=0$ in $L^2$; hence $S=0$.
The coefficient of $\mathbf K_t^{\mathbf E}$ at this element is
$\langle P_t^{\mathbf E}\boldsymbol\phi_j,
\boldsymbol\phi_i\rangle_{L^2}=e^{-t\lambda_j}\delta_{ij}$.
This proves the second series in $L^2(M\times M)$, and at $t=\displaystyle\frac{\tau}{2}$ gives
\[
 \sum_{j=1}^{\infty}e^{-\tau\lambda_j}
       =\|\mathbf K_{\frac{\tau}{2}}^{\mathbf E}\|_{L^2(M\times M)}^2<\infty.
\]
The connection Laplacian of the product bundle is the sum of the
Laplacians in the two variables; denote it by $L_{\mathrm{prod}}$. Its eigenvalue on the diagonal
term is $2\lambda_j$. For $a,b\geq0$, the series tails after
applying $(I+L_{\mathrm{prod}})^a\partial_t^b$ have squared norm
bounded, uniformly for $t\geq\tau$, by
\[
 \sum_{j>N}(1+2\lambda_j)^{2a}\lambda_j^{2b}e^{-2\tau\lambda_j}
 \leq C'_{a,b,\tau}\sum_{j>N}e^{-\tau\lambda_j}
 \longrightarrow0.
\]
Elliptic regularity and Sobolev embeddings on the closed manifold
$M\times M$ prove uniform convergence of all kernel derivatives.
Finally, on the orthogonal complement of the kernel,
$e^{-t\lambda_j}\leq e^{-t\lambda_+}$; Parseval gives the convergence
estimate.
\end{proof}

\begin{proposition}[Dirichlet and Neumann realizations]
\label{prop:espectro-compacto-frontera-formas}
Now let $M$ be a compact, connected Riemannian manifold of positive
dimension with nonempty smooth boundary. For scalar functions, the form
\[
 \mathfrak q(u,v):=\int_M\langle du,dv\rangle_{\mathbf g}\,
                              d\lambda_{\mathbf g}
\]
with domain $V_D:=W_0^{1,2}(M)$, respectively $V_N:=W^{1,2}(M)$,
defines the nonnegative self-adjoint realizations $L_D$ and $L_N$.
Both have compact resolvent, an orthonormal eigenfunction basis,
eigenvalues of finite multiplicity tending to $+\infty$, and the
series domain and minimax formulas of this section.
Eigenfunctions are smooth in the interior.
We have $\lambda_1^D>0$, while $\lambda_1^N=0$ is simple.

The following distributional equations are interpreted on $\operatorname{Int}M$. Boundary conditions are encoded by the form domain and the variational identity.

Their domains are characterized by
\[
 \begin{split}
 \mathcal D(L_D)&=\{u\in W_0^{1,2}(M)\mid
                        -\Delta_{\mathbf g}u\in L^2(M)\},\\
 \mathcal D(L_N)&=\left\{u\in W^{1,2}(M)\ \middle|\
 \begin{array}{l}
 \text{there exists }f\in L^2(M)\text{ such that}\\
 \mathfrak q(u,v)=\langle f,v\rangle_{L^2}
                       \text{ for every }v\in W^{1,2}(M)
 \end{array}\right\}.
 \end{split}
\]
In the second case, $L_Nu=f=-\Delta_{\mathbf g}u$ in $\mathcal D'(\operatorname{Int}M)$.
For functions smooth up to the boundary, these characterizations
are equivalent to $u|_{\partial M}=0$ and $\partial_\nu u=0$,
respectively, where $\nu$ is the outward normal.
\end{proposition}

\begin{proof}
The spaces $V_D$ and $V_N$ are dense in $L^2(M)$ and complete for
$\|u\|_{L^2}^2+\mathfrak q(u,u)$. The form representation theorem
constructs $L_D$ and $L_N$. Both form domains embed compactly
into $L^2(M)$ by Rellich with boundary.
Lemma~\ref{lem:criterio-forma-compacta-eigenbase} and
Proposition~\ref{prop:rayleigh-minimax-forma-compacta} give the
spectral assertions. The distributional eigenfunction equation
and local elliptic regularity give interior smoothness.

The domain characterization for a form-associated operator is
$\mathfrak q(u,v)=\langle f,v\rangle$ for every $v$ in the form
domain. In the Dirichlet case, testing first against
$C_c^\infty(\operatorname{Int}M)$ is equivalent to
$-\Delta_{\mathbf g}u=f$; density of that space in $W_0^{1,2}(M)$
extends the identity to all of $V_D$. In the Neumann case, the
characterization requires every element of $W^{1,2}(M)$, and
interior tests identify $f$ with $-\Delta_{\mathbf g}u$.
For smooth $u$, Green's formula turns the Neumann identity into
$\int_{\partial M}\partial_\nu u\,\overline v\,d\lambda_{\partial M}=0$
for every smooth $v$. Every smooth boundary test extends smoothly
using a collar and a partition of unity, so this is equivalent
to $\partial_\nu u=0$. The Dirichlet condition is zero trace,
by Theorem~\ref{teo:caracterizacion-W0-trazas}.

A kernel element of either realization has zero weak gradient
and is constant because $M$ is connected. In $V_D$, the trace
of that constant must vanish on the nonempty boundary, so the
kernel is zero. In $V_N$, the kernel is precisely the constants.
Discreteness of the spectrum gives the assertions about the
first eigenvalues.
\end{proof}

\begin{example}[Circle and interval]
\label{ej:espectro-circulo-intervalo}
On the circle of length $\ell$, parametrized by arc length,
the normalized complex eigenvectors are
$\ell^{-\displaystyle\frac{1}{2}}e^{2\pi \displaystyle\frac{ikx}{\ell}}$, $k\in\mathbb Z$, with eigenvalues
$(2\pi \displaystyle\frac{k}{\ell})^2$. This follows by differentiating twice and using
the periodic Fourier basis. Zero is simple, and each positive
eigenvalue has multiplicity two.

On $[0,\ell]$, the equation $-\phi''=\lambda\phi$ and Dirichlet
conditions give $\lambda_k^D=(\pi \displaystyle\frac{k}{\ell})^2$, $k\geq1$, with
eigenfunctions $(\displaystyle\frac{2}{\ell})^{\displaystyle\frac{1}{2}}\sin(\pi \displaystyle\frac{kx}{\ell})$.
For Neumann conditions, the constant $\ell^{-\displaystyle\frac{1}{2}}$ has eigenvalue
zero, and $(\displaystyle\frac{2}{\ell})^{\displaystyle\frac{1}{2}}\cos(\pi \displaystyle\frac{kx}{\ell})$ has eigenvalue
$(\pi \displaystyle\frac{k}{\ell})^2$, $k\geq1$.
The sine and cosine bases follow from odd and even extensions
to a periodic interval. These extensions and Fourier Parseval
prove completeness of the bases.
\end{example}

\chapter{Sobolev spaces on manifolds of bounded geometry}
\label{cap:geometria-acotada}

An estimate obtained in a chart may be correct and yet provide no global information: its constant may grow as the center of the chart tends to infinity. On a compact manifold, choosing a finite cover suffices. On a noncompact manifold, we need another property: the ability to repeat the same local argument with radii and constants independent of the point.

Bounded geometry makes this possible. The injectivity radius provides normal charts of a common size; bounds on the curvature and its derivatives control the metric in these charts. From these data we will construct covers of uniformly bounded multiplicity and partitions of unity whose derivatives remain bounded. These are the conditions that allow us to sum the local estimates without losing control of the norm.

We will then compare spaces defined by localization with those defined using covariant derivatives or powers of the Laplacian. In this comparison, we must distinguish integer regularity from fractional regularity and, within the latter, the different orders of summation defining Besov and Triebel--Lizorkin spaces. We will retain these distinctions in interpolation and duality arguments as well.

The manifolds in this chapter are assumed to have no boundary. The case with boundary is developed in Chapter~\ref{cap:trazas-geometria-acotada}, where normal charts near the boundary are replaced by Fermi charts, and Euclidean norms by restriction norms on the half-space.

\begin{semblanzaHistorica}{Uniformity in noncompact spaces}
On a compact manifold, a finite cover allows the distortions of the charts and the metric to be absorbed into a single constant. This advantage disappears in the noncompact case. Bounded geometry replaces finiteness by uniformity: charts of controlled size, bounded multiplicity, and estimates independent of the center. Eichhorn's work and the framework of uniformly regular manifolds developed by Amann provide two closely related ways of turning this geometric uniformity into a stable theory of function spaces~\cite{Eichhorn2007,Amann2025FunctionSpaces}.
\end{semblanzaHistorica}

\section{Riemannian manifolds of bounded geometry}
In the analytic study of a Riemannian manifold, it is essential to have local coordinate systems in which the geometry is controlled and comparable to its Euclidean counterpart. One of the most effective tools for this purpose is normal coordinates, constructed using the Riemannian exponential map at a point. These coordinates give a particularly simple representation of the metric near their center and facilitate the transfer of classical techniques from $\mathbb{R}^n$ to manifolds.

However, this procedure has limitations: at some distance, the exponential map may cease to be a diffeomorphism onto its image because conjugate points occur or distinct geodesics issuing from the same point intersect. This leads to the need to quantify the largest radius within which normal coordinates behave appropriately.

This parameter is the \emph{injectivity radius} of the manifold, which plays a central role in geometry and analysis. A positive lower bound for this radius allows us to use normal charts of a common size. To control their metric coefficients and derivatives as well, we will need the curvature bounds introduced below. Both conditions enter into comparisons of Sobolev, Besov, and Triebel--Lizorkin norms.

Recall Definition~\ref{def: radio de inyectividad}: $\operatorname{inj}(p)$ measures the largest radius of a tangent ball on which $\exp_p$ is a diffeomorphism onto its image, and $\operatorname{inj}(M)=\displaystyle\inf_{p\in M}\operatorname{inj}(p)$. Pointwise positivity does not imply a uniform positive lower bound on a noncompact manifold.
\begin{remark}\label{obs: geometria acotada implica completa}
Let $M$ be a Riemannian manifold without boundary. If $\operatorname{inj}(M)>0$, then every connected component of $M$ is complete. Indeed, choose $0<r<\operatorname{inj}(M)$. Every unit-speed geodesic defined on a maximal interval $[0,b)$, with $b<\infty$, can be extended from any time $t_0<b$ for at least a time $r$, since $\exp_{\gamma(t_0)}$ is defined on the ball of radius $r$. Taking $b-r<t_0<b$ gives an extension beyond $b$, a contradiction. Thus, each component is geodesically complete; the Hopf--Rinow Theorem~\ref{teo:variedades-riemannianas-hopf-rinow} implies its metric completeness.
\end{remark}
Since $\operatorname{inj}(M)\geq 0$, there are only two possibilities: $\operatorname{inj}(M)=0$ or $\operatorname{inj}(M)>0$. The first behavior can occur, for example, on manifolds containing geodesic loops of arbitrarily small length, as illustrated in the following figure.

\begin{figura}[htbp]
    \includegraphics[scale=1]{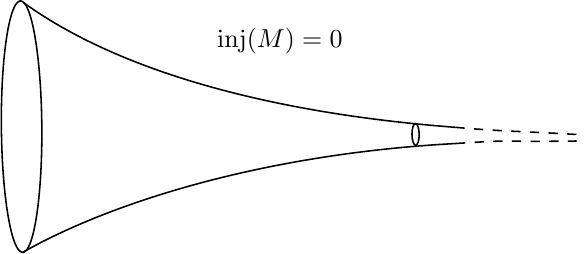}
    \caption{Schematic representation of a manifold with geodesic loops of arbitrarily small length, a phenomenon that can give rise to zero injectivity radius.}
    \label{fig:radio-inyectividad-cero}
\end{figura}
\begin{remark}\label{obs:geometria-acotada-variedad-riemanniana-su-radio-inyectividad-satisface}
Let $(M,\mathbf{g})$ be a Riemannian manifold without boundary, and suppose that its injectivity radius satisfies $\operatorname{inj}(M)>0$. Then, for each $p\in M$ and every $0<r<\operatorname{inj}(M)$, the exponential map
\[
\operatorname{exp}_p\colon B_{\mathbf{g}_p}(0,r)\subseteq T_pM \longrightarrow B_{\mathbf{g}}(p,r)\subseteq M
\]
is a diffeomorphism.

In particular, every point $q\in B_{\mathbf{g}}(p,r)$ can be joined to $p$ by a unique minimizing geodesic of length less than $r$, given by
\[
\gamma(t)=\operatorname{exp}_p(t v), \qquad t\in[0,1],
\]
where $v\in T_pM$ satisfies $\|v\|_{\mathbf{g}_p}<r$ and $\operatorname{exp}_p(v)=q$.
\end{remark}
To quantify how closely the Riemannian metric approximates the Euclidean metric in normal coordinates, we will study the Taylor expansion of its coefficients. The following proposition makes this comparison explicit.
\begin{proposition}\label{prop: expansion de taylor de la metrica en coordenadas normales}
Let $(M,\mathbf{g})$ be a Riemannian manifold with or without boundary, and let $p\in\operatorname{Int}(M)$. If $(U,\phi)$ is a normal coordinate chart centered at $p$, with $\phi(U)$ convex and $\phi(p)=0$, then, for every $q\in U$, writing $x=(x^{1},\dots,x^{n})=\phi(q)$, we have the following Taylor expansion:
\[
(g_{ij}\circ\phi^{-1})(x)
=\delta_{ij}
-\frac{1}{3}R_{iklj}(p)x^{k}x^{l}
+R_{2,0,g_{ij}\circ\phi^{-1}}(x),
\]
where, for every $\varepsilon>0$, there exists $\delta_\varepsilon>0$ such that
\[
 |R_{2,0,g_{ij}\circ\phi^{-1}}(x)|
 \leq\varepsilon\|x\|^2
 \quad\text{if }\|x\|<\delta_\varepsilon.
\]
\end{proposition}

\begin{proof}
Since \(p\in\operatorname{Int}(M)\), we may consider the manifold without boundary \(\operatorname{Int}(M)\) with the restricted metric. By Proposition~\ref{propiedades de las coordenadas normales} applied to this manifold, in normal coordinates centered at $p$ we have
\[
g_{ij}(p)=\delta_{ij},\qquad \Gamma^{k}_{ij}(p)=0.
\]
By Corollary~\ref{compatibilidad levi civita metrica en coordenadas} (metric compatibility),
\begin{equation}\label{compatibiilidad metrica derivadas}\partial_{k}g_{ij}=g_{mj}\Gamma^{m}_{ik}+g_{im}\Gamma^{m}_{jk}.
\end{equation}
Evaluating at $p$ and using $g_{mj}(p)=\delta_{mj}$, $g_{im}(p)=\delta_{im}$, and the vanishing of the Christoffel symbols at $p$, we obtain
\[
(\partial_{k}g_{ij})(p)=0.
\]
Thus, the first-order derivatives vanish at $p$, and the corresponding linear term does not appear in the Taylor expansion.

Now, for arbitrary $v=v^{i}\boldsymbol{\partial}_{i}|_{p}\in T_{p}M$, Proposition~\ref{propiedades de las coordenadas normales} gives $\gamma_{v}(t)=(tv^{1},\dots,tv^{n})$, so the geodesic equation
\[\frac{d^{2}x^{k}}{dt^{2}}(t)+\frac{dx^{i}}{dt}(t)\frac{dx^{j}}{dt}(t)\Gamma_{ij}^{k}(x(t))=0\]
with $x(t)=\gamma_{v}(t)$ becomes
\[\Gamma_{ij}^{k}(tv)v^{i}v^{j}=0.\]

Differentiating with respect to $t$ gives:
\[0=\frac{d}{dt}(\Gamma_{ij}^{k}(tv)v^{i}v^{j})=\left(\frac{d}{dt}(\Gamma_{ij}^{k}(tv))\right) v^{i}v^{j}=\partial_{l}(\Gamma_{ij}^{k})(tv)\frac{d(tv^{l})}{dt}v^{i}v^{j}=\partial_{l}(\Gamma_{ij}^{k})(tv)v^{i}v^{j}v^{l}.\]

Evaluating at $t=0$ gives
\[\partial_{l}\Gamma_{ij}^{k}(p)v^{i}v^{j}v^{l}=0.\]
This is a homogeneous polynomial of degree $3$ that vanishes identically (since $v$ was arbitrary). The coefficient of each term $v^{i}v^{j}v^{l}$ is the sum over $6$ permutations of these indices ($i,j$ and $l$). Since $\partial_{l}\Gamma_{ij}^{k}=\partial_{l}\Gamma_{ji}^{k}$, we obtain
\begin{equation}\label{parciales christoffel coordenadas normales en p}
 \partial_{l}\Gamma_{ij}^{k}(p)+\partial_{i}\Gamma_{jl}^{k}(p)+\partial_{j}\Gamma_{li}^{k}(p)=0.\end{equation}

We now differentiate (\ref{compatibiilidad metrica derivadas}) once more:
\[
\partial_{l}\partial_{k}g_{ij}
=\partial_{l}\big(g_{mj}\Gamma^{m}_{ik}+g_{im}\Gamma^{m}_{jk}\big)=
(\partial_{l}g_{mj})\Gamma^{m}_{ik}
+g_{mj}(\partial_{l}\Gamma^{m}_{ik})
+(\partial_{l}g_{im})\Gamma^{m}_{jk}
+g_{im}(\partial_{l}\Gamma^{m}_{jk}).
\]
Evaluating at $p$ and again using the vanishing of the Christoffel symbols and the partial derivatives of order $1$ at $p$, only terms containing derivatives of the Christoffel symbols remain. Consequently,
\[
(\partial_{l}\partial_{k}g_{ij})(p)
=g_{mj}(p)(\partial_{l}\Gamma^{m}_{ik})(p)
+g_{im}(p)(\partial_{l}\Gamma^{m}_{jk})(p).
\]
Finally, since the metric coefficients evaluated at $p$ are Kronecker deltas, we are left with:
\begin{equation}\label{segunda derivada gij en terminos de christoffel}
(\partial_{l}\partial_{k}g_{ij})(p)
=(\partial_{l}\Gamma^{j}_{ik})(p)
+(\partial_{l}\Gamma^{i}_{jk})(p).
\end{equation}

On the other hand, Proposition~\ref{1,3 curvatura en coordenadas} gives
\[R_{ijk}{}^{l}=\partial_i\Gamma_{jk}^{l}-\partial_{j}\Gamma_{ik}^{l}+\Gamma_{jk}^{m}\Gamma_{im}^{l}-\Gamma_{ik}^{m}\Gamma_{jm}^{l},\]

so evaluating at $p$ and using the vanishing of the Christoffel symbols at $p$, we obtain

\[
R_{ikl}{}^{j}(p)
= \partial_{i}\Gamma^{j}_{kl}(p)
- \partial_{k}\Gamma^{j}_{il}(p),
\qquad
R_{ilk}{}^{j}(p)
= \partial_{i}\Gamma^{j}_{lk}(p)
- \partial_{l}\Gamma^{j}_{ik}(p).
\]
Adding these terms gives
\[R_{ikl}{}^{j}(p)+R_{ilk}{}^{j}(p)=2\partial_{i}\Gamma_{kl}^{j}(p)-\partial_{l}\Gamma_{ik}^{j}(p)-\partial_{k}\Gamma_{il}^{j}(p)\]
and using (\ref{parciales christoffel coordenadas normales en p}) with $k$ and $j$ interchanged, we obtain
\[ \partial_{l}\Gamma_{ik}^{j}(p)+\partial_{i}\Gamma_{kl}^{j}(p)+\partial_{k}\Gamma_{li}^{j}(p)=0,\]
so that
\[R_{ikl}{}^{j}+R_{ilk}{}^{j}=3\partial_{i}\Gamma_{kl}^{j}(p).\]

Since $g_{ij}=\delta_{ij}$ at $p$, the preceding identity can be written in terms of the Riemann curvature tensor as follows:
\[R_{iklj}(p)+R_{ilkj}(p)=R_{ikl}{}^{j}(p)+R_{ilk}{}^{j}(p)=3\partial_{i}\Gamma_{kl}^{j}(p).\]
Interchanging $i$ and $l$ gives
\[
R_{likj}(p) + R_{lkij}(p) = 3 \partial_{l}\Gamma^{ j}_{ik}(p)\]
and interchanging $i$ and $l$, and then $i$ and $j$, gives
\[R_{ljki}(p) + R_{lkji}(p) = 3 \partial_{l}\Gamma^{ i}_{jk}(p).\]

Using identity (\ref{segunda derivada gij en terminos de christoffel}), we obtain
\[3\partial_{l}\partial_{k}g_{ij}(p)=R_{likj}(p) + R_{lkij}(p) + R_{ljki}(p) + R_{lkji}(p).\]
By Corollary~\ref{simetrías tensor curvatura coordenadas}, $R_{lkij}(p)=-R_{lkji}(p)$ and $R_{ljki}(p)=R_{kilj}(p)$, so
\[\partial_{l}\partial_{k}g_{ij}(p)=\frac{1}{3}\left(R_{likj}(p)+R_{kilj}(p)\right)=-\frac{1}{3}(R_{iklj}(p)+R_{ilkj}(p)),\]
where the last equality uses $R_{likj}(p)=-R_{ilkj}(p)$ and $R_{kilj}(p)=-R_{iklj}(p)$.

By (Taylor's) Theorem~\ref{taylor multivariable multiindices},
\[
\begin{aligned}
(g_{ij}\circ\phi^{-1})(x)
={}&\delta_{ij}+\partial_kg_{ij}(p)x^k
+\frac12\partial_k\partial_lg_{ij}(p)x^kx^l
+R_{2,0,g_{ij}\circ\phi^{-1}}(x)\\
={}&\delta_{ij}
-\frac16\bigl(R_{iklj}(p)+R_{ilkj}(p)\bigr)x^kx^l
+R_{2,0,g_{ij}\circ\phi^{-1}}(x)\\
={}&\delta_{ij}-\frac13R_{iklj}(p)x^kx^l
+R_{2,0,g_{ij}\circ\phi^{-1}}(x).
\end{aligned}
\]
In the last equality, the dummy indices $k$ and $l$ are interchanged in the second summand.
\end{proof}

Thus, curvature determines the first deviation of the metric from the Euclidean metric in normal coordinates. Its covariant derivatives enter into higher-order expansions; see \cite{MullerSchubertVanDeVen1999}. However, an expansion at each center does not by itself provide a uniform radius or bound for its remainders. To obtain this control, we will work directly with the Jacobi equation. The hypotheses that make this possible are as follows.
\begin{definition}[Manifolds of bounded geometry]\label{def: geometria acotada}\index{manifold!of bounded geometry}\index{bounded geometry@bounded geometry}\glsadd{geometria-acotada}
Let $(M,\mathbf{g})$ be a Riemannian manifold without boundary. We say that $M$ has bounded geometry if the following two conditions hold:
\begin{enumerate}
\item $\operatorname{inj}(M)>0$.
\item For each $k\in\mathbb N_0$, there exists $C_k>0$ such that
\[
\displaystyle\sup_{x\in M}|\nabla^k\mathbf{R}(x)|_{\mathbf{g}}\leq C_k.
\]
Here $\mathbf{R}\in\Gamma(T^{(1,3)}(TM))$ is the curvature tensor of Definition~\ref{def:variedades-riemannianas-variedad-5}, $\nabla^k\mathbf{R}$ is its $k$th derivative for the tensor connection induced by the Levi--Civita connection, and $|\cdot|_{\mathbf{g}}$ is the pointwise norm induced by $\mathbf{g}$ and its dual metric on all tensor factors.
\end{enumerate}
\end{definition}

\begin{definition}[Bounded geometry of finite order]
\label{def:geometria-acotada-orden-k}
\label{def:condiciones-I-Bk-eichhorn}
\index{bounded geometry!of finite order}
Let $(M,\mathbf g)$ be a smooth Riemannian manifold without boundary,
and let $k\in\mathbb N_0$. It has \emph{bounded geometry of order
$k$} if the following conditions hold:
\[
 \begin{aligned}
 (I)\qquad &\operatorname{inj}(M,\mathbf g)\geq i_0
                     &&\text{for some constant }i_0>0,\\
 (B_k)\qquad &\sup_{x\in M}|(\nabla^{\mathbf g})^j\mathbf R(x)|_{\mathbf g}
                   \leq K_j
                     &&\text{for }0\leq j\leq k,
 \end{aligned}
\]
with constants $K_0,\ldots,K_k>0$. Derivatives of $\mathbf R$ are
taken in the corresponding tensor bundle, as in
Definition~\ref{def: geometria acotada}. The integer $k$ counts the
curvature derivatives being controlled; the metric remains smooth.
\end{definition}

Bounded geometry, without an order specified, means bounded geometry
of all orders. It is equivalent to satisfying $(I)$ and $(B_k)$ for
every $k\in\mathbb N_0$; the constants $K_j$ may depend on $j$.
Condition $(I)$ implies completeness of each component by
Remark~\ref{obs: geometria acotada implica completa}, even when only
$(B_0)$ is required. Thus completeness is not an additional hypothesis
in these definitions for manifolds without boundary. Conditions
$(I)$ and $(B_k)$ agree with those used by Eichhorn; see
\cite[\S 2]{EichhornSlice}.

\begin{remark}[Euclidean open subsets and manifold boundary]
\label{obs:abierto-euclidiano-no-geometria-acotada}
An open subset of $\mathbb R^n$ equipped with the restricted Euclidean metric need not have bounded geometry. For example, let $M=B_{\mathrm{euc}}(0,1)$, regarded as a manifold without boundary. Its curvature is zero, but, for $x\in M$, the exponential map is $\exp_x(v)=x+v$ only as long as the corresponding segment remains in the ball. Consequently,
\[
 \operatorname{inj}(x)=1-|x|,
 \qquad \operatorname{inj}(M)=0.
\]
The radial geodesic moving toward the sphere reaches it in finite time and cannot be extended within $M$. Completeness does not appear as a third hypothesis in the preceding definition: it follows from the positive injectivity radius by Remark~\ref{obs: geometria acotada implica completa}, and this is precisely the condition that fails here.

The sphere is the boundary of this open set as a Euclidean subset; it is not a boundary included in the manifold $M$. In contrast, $\overline{B}_{\mathrm{euc}}(0,1)$, with its smooth structure as a compact manifold with boundary and the same Euclidean metric, belongs to the theory with boundary of Chapter~\ref{cap:trazas-geometria-acotada}. These are two distinct geometric objects. The metric matters as well: replacing it by a complete metric can change the conclusion concerning completeness. Thus, the Sobolev and Rellich theorems for arbitrary Euclidean open subsets, particularly the zero-extension arguments for $W_0^{m,p}(M)$, do not follow simply by applying the global bounded-geometry theorems.
\end{remark}

Since throughout this chapter we will work with \textit{\textbf{normal coordinates}} (also called \textit{\textbf{geodesic coordinates}}), as well as with coordinates in general, we introduce some terminology that will allow us to refer to them more clearly.

\begin{definition}[Geodesic atlas]\label{def: atlas geodesico}\index{geodesic atlas@geodesic atlas}
Let $(M,\mathbf{g})$ be a Riemannian manifold without boundary.

We will denote by \[
\mathcal{A}^{\operatorname{geo}}=(U_{i}^{\operatorname{geo}},\phi_{i}^{\operatorname{geo}})_{i\in J}
\] an atlas consisting of normal coordinate charts (see Proposition~\ref{propiedades de las coordenadas normales}).

More precisely, for each $i\in J$, there exist a point $p_{i}\in M$, a radius $r_{i}>0$, and a linear isomorphism \[
B_{p_{i}}\colon \mathbb{R}^{n}\longrightarrow T_{p_{i}}M,
\] taking the standard basis of $\mathbb{R}^{n}$ to an orthonormal basis of $T_{p_{i}}M$, such that
\[
U_{i}^{\operatorname{geo}}=B_{\mathbf{g}}(p_i,r_i),
\qquad
\phi_{i}^{\operatorname{geo}}
=
B_{p_{i}}^{-1}\circ \exp_{p_{i}}^{-1}.
\]

We say that $\mathcal{A}^{\operatorname{geo}}$ is a \textit{uniform geodesic atlas} if there exists $r>0$ such that $r<\operatorname{inj}(M)$ and \[
U_{i}^{\operatorname{geo}}=B_{\mathbf{g}}(p_i,r)
\] for every $i\in J$.

\end{definition}

The condition that the covariant derivatives of the curvature are bounded can be expressed in other ways which, together with $\operatorname{inj}(M)>0$, yield useful equivalent characterizations of manifolds of bounded geometry.
\begin{theorem}[Equivalent characterizations of bounded geometry]
\label{teo: equivalencias de geometria acotada}\index{bounded geometry@bounded geometry!equivalent characterizations}
Let $(M,\mathbf{g})$ be a Riemannian manifold without boundary such that $\operatorname{inj}(M)>0$. The following conditions are equivalent:
\begin{enumerate}
\item For each $k\in \mathbb{N}_{0}$, there exists $C_{k}>0$ such that
\[
 |\nabla^{k}\mathbf{R}|_{\mathbf{g}}\leq C_{k}.
 \]
\item There exists \(r_0\), with \(0<r_0<\operatorname{inj}(M)\), such that the full family of normal charts of radius \(r_0\) satisfies the following condition: for each \(k\in\mathbb N_0\), there exists \(C_{k,r_0}>0\) such that, in every chart, the components of $\mathbf{g}$ and of its inverse satisfy
\[
 |D^{\alpha}g_{ij}|\leq C_{k,r_0},
 \qquad
 |D^{\alpha}g^{ij}|\leq C_{k,r_0},
 \qquad |\alpha|\leq k.
 \]
By restriction, the same estimates hold for every radius \(0<r\leq r_0\).
\end{enumerate}

Whenever either condition holds, transitions between normal charts of any fixed radius \(0<r\leq r_0\) have uniformly bounded jets of every order: for each \(k\in\mathbb N_0\), there exists \(A_{k,r}>0\) such that, if the charts \(i\) and \(j\) overlap,
\[
\bigl|D^\alpha(\phi_j^{\operatorname{geo}}\circ
(\phi_i^{\operatorname{geo}})^{-1})(x)\bigr|
\leq A_{k,r},
\qquad |\alpha|\leq k,
\]
throughout the domain of the transition. The same assertion holds with \(i\) and \(j\) interchanged.
\end{theorem}
\begin{proof}
Fix a common radius
\[
0<\rho<\min\{\operatorname{inj}(M),1\}.
\]
When necessary, we will decrease it once; all normal charts of radius \(\rho\) will still cover \(M\). For each \(p\in M\), choose a linear isometry \(B_p\colon \mathbb R^n\longrightarrow T_pM\) and write
\[
\Psi_p(v):=\exp_p(B_pv).
\]
We denote by \(g_{ij}^{(p)}\) and \(g_{(p)}^{ij}\) the components of \(\Psi_p^*\mathbf g\) and of its inverse metric, respectively, in the standard basis of \(\mathbb R^n\).

First assume (1). For \(v\in B(0,\rho)\) and \(i\in\{1,\dots,n\}\), the variation
\[
\alpha_i(s,t):=\exp_p\bigl(tB_p(v+se_i)\bigr)
\]
produces the Jacobi field
\(\mathbf{J}_i(t)=\partial_s\alpha_i(0,t)\)
along \(\gamma_v(t)=\Psi_p(tv)\). Its initial data and equation are
\[
\mathbf{J}_i(0)=0,
\qquad
\frac{D\mathbf{J}_i}{dt}(0)=B_pe_i,
\qquad
\frac{D^2\mathbf{J}_i}{dt^2}
+\mathbf{R}(\mathbf{J}_i,\dot\gamma_v)\dot\gamma_v=0.
\]
Moreover, \(\mathbf{J}_i(1)=d\Psi_p(v)e_i\). Parallel transport the basis \((B_pe_a)_{a=1}^n\) along \(\gamma_v\). In this frame, the equation is a linear system
\[
y_i''(t)+\mathcal R_{p,v}(t)y_i(t)=0,
\qquad y_i(0)=0,
\qquad y_i'(0)=e_i,
\]
and \(\|\mathcal R_{p,v}(t)\|\leq C_0\|v\|^2\). Writing it as a first-order system and applying Lemma~\ref{lema: gronwall} gives a bound for \(\|y_i\|+\|y_i'\|\) independent of \(p\) and \(v\). The integral form of the equation also gives
\[
\|y_i(t)-te_i\|
\leq
C_0\|v\|^2\int_0^t(t-\tau)\|y_i(\tau)\|\,d\tau.
\]
Decrease \(\rho\), depending only on \(C_0\) and \(n\), so that the right-hand side at \(t=1\) is less than \(\frac{1}{2\sqrt n}\) for every \(i\). If \(\mathbf{Y}(v)\) is the matrix whose columns are \(y_i(1)\), then
\[
\|\mathbf{Y}(v)-I\|_{\mathrm{op}}
\leq
\|\mathbf{Y}(v)-I\|_{\mathrm{HS}}
<
\frac12.
\]
By the Neumann series, \(\mathbf{Y}(v)\) is invertible and \(\|\mathbf{Y}(v)^{-1}\|_{\mathrm{op}}\leq2\), uniformly in \(p\) and \(v\). Since the parallel frame is orthonormal,
\[
g_{ij}^{(p)}(v)=\langle \mathbf{J}_i(1),\mathbf{J}_j(1)\rangle_{\mathbf{g}},
\]
which proves the uniform ellipticity of the coefficient matrix \((g_{ij}^{(p)})\).

We now justify the derivatives of arbitrary order. Fix \(m\in\mathbb N_0\) and vectors \(w_0,\ldots,w_m\in\{e_1,\ldots,e_n\}\) and, for sufficiently small \(\boldsymbol s=(s_0,\ldots,s_m)\), consider the multiparameter variation
\[
 \alpha(\boldsymbol s,t)
 :=\exp_p\!\left(tB_p\left(v+\displaystyle\sum_{a=0}^m s_aw_a\right)\right),
 \qquad 0\leq t\leq1.
\]
Let \(\boldsymbol{\mathcal{E}}:=\alpha^*TM\), and denote by \(\mathbb D\) the connection on \(\boldsymbol{\mathcal{E}}\) induced by the Levi--Civita connection. Set
\[
 D_t:=\mathbb D_{\boldsymbol{\partial}_t},\qquad
 D_a:=\mathbb D_{\boldsymbol{\partial}_{s_a}},\qquad
 \mathbf{T}:=d\alpha(\boldsymbol{\partial}_t),\qquad
 \mathbf{A}_a:=d\alpha(\boldsymbol{\partial}_{s_a}).
\]
For each \(\boldsymbol s\), the curve \(t\mapsto\alpha(\boldsymbol s,t)\) is geodesic; moreover, since the Levi--Civita connection is torsion-free,
\[
 D_t\mathbf{T}=0,
 \qquad
 D_a\mathbf{T}=D_t\mathbf{A}_a.
\]
With the convention \(\mathbf{R}(\mathbf{X},\mathbf{Y})\mathbf{Z}=\nabla_{\mathbf{X}}\nabla_{\mathbf{Y}}\mathbf{Z}-\nabla_{\mathbf{Y}}\nabla_{\mathbf{X}}\mathbf{Z}-\nabla_{[\mathbf{X},\mathbf{Y}]}\mathbf{Z}\), the commutators on the pullback bundle are precisely
\[
 (D_aD_t-D_tD_a)\mathbf{Z}=\mathbf{R}(\mathbf{A}_a,\mathbf{T})\mathbf{Z},
 \qquad
 (D_aD_b-D_bD_a)\mathbf{Z}=\mathbf{R}(\mathbf{A}_a,\mathbf{A}_b)\mathbf{Z}.
\]
In particular, differentiating \(D_t\mathbf{T}=0\) in direction \(s_a\) gives, with this sign,
\[
 D_t^2\mathbf{A}_a+\mathbf{R}(\mathbf{A}_a,\mathbf{T})\mathbf{T}=0.
\]

We will carry out the induction without permuting covariant derivatives. For an ordered word \(I=(a_1,\ldots,a_k)\), write
\[
 |I|=k,
 \qquad
 D_I:=D_{a_k}\cdots D_{a_1},
 \qquad
 I^-:=(a_1,\ldots,a_{k-1}),
\]
and set \(D_\varnothing=\operatorname{id}\). Introduce the operator
\[
 \mathcal L\mathbf{Z}:=D_t^2\mathbf{Z}+\mathbf{R}(\mathbf{Z},\mathbf{T})\mathbf{T},
 \qquad
 \mathcal L\colon\Gamma(\boldsymbol{\mathcal{E}})\longrightarrow\Gamma(\boldsymbol{\mathcal{E}}).
\]
A direct calculation using the preceding commutators gives, for every \(\mathbf{Z}\in\Gamma(\boldsymbol{\mathcal{E}})\),
\[
 \mathcal L(D_a\mathbf{Z})=D_a(\mathcal L\mathbf{Z})+\mathcal C(\mathbf{A}_a,\mathbf{Z}),
\]
where \(\mathcal C(\mathbf{A},\mathbf{Z})\in\Gamma(\boldsymbol{\mathcal{E}})\) is
\begin{align*}
 \mathcal C(\mathbf{A},\mathbf{Z})={}&
 -(\nabla_{\mathbf{T}}\mathbf{R})(\mathbf{A},\mathbf{T})\mathbf{Z}-\mathbf{R}(D_t\mathbf{A},\mathbf{T})\mathbf{Z}-2\mathbf{R}(\mathbf{A},\mathbf{T})D_t\mathbf{Z}\\
 &-(\nabla_{\mathbf{A}}\mathbf{R})(\mathbf{Z},\mathbf{T})\mathbf{T}-\mathbf{R}(\mathbf{Z},D_t\mathbf{A})\mathbf{T}-\mathbf{R}(\mathbf{Z},\mathbf{T})D_t\mathbf{A}.
\end{align*}
Indeed,
\[
 D_t^2D_a\mathbf{Z}
 =D_aD_t^2\mathbf{Z}-(\nabla_{\mathbf{T}}\mathbf{R})(\mathbf{A}_a,\mathbf{T})\mathbf{Z}-\mathbf{R}(D_t\mathbf{A}_a,\mathbf{T})\mathbf{Z}
   -2\mathbf{R}(\mathbf{A}_a,\mathbf{T})D_t\mathbf{Z},
\]
and the derivative of \(\mathbf{R}(\mathbf{Z},\mathbf{T})\mathbf{T}\), using \(D_a\mathbf{T}=D_t\mathbf{A}_a\), contributes precisely the remaining three terms. This verifies both the signs and the tensor spaces: each summand of \(\mathcal C(\mathbf{A},\mathbf{Z})\) is a section of \(\boldsymbol{\mathcal{E}}\).

For a distinguished letter \(b\in\{0,\ldots,m\}\), define
\[
 \mathbf{Z}_{I;b}:=D_I\mathbf{A}_b\in\Gamma(\boldsymbol{\mathcal{E}}),
 \qquad
 \widetilde{\mathbf{Q}}_{\varnothing;b}:=0,
\]
and, if \(I=(I^-,a)\), define recursively
\[
 \widetilde{\mathbf{Q}}_{I;b}
 :=D_a\widetilde{\mathbf{Q}}_{I^-;b}
   +\mathcal C(\mathbf{A}_a,\mathbf{Z}_{I^-;b})
 \in\Gamma(\boldsymbol{\mathcal{E}}).
\]
We prove the recurrence by induction on $|I|$. For $I=\varnothing$, the Jacobi equation gives
\[
\mathcal L\mathbf Z_{\varnothing;b}
=\mathcal L\mathbf A_b=0
=\widetilde{\mathbf Q}_{\varnothing;b}.
\]
Now suppose that $\mathcal L\mathbf Z_{I^-;b}=\widetilde{\mathbf Q}_{I^-;b}$ and write $I=(I^-,a)$. Since $\mathbf Z_{I;b}=D_a\mathbf Z_{I^-;b}$, the preceding commutation identity gives, term by term,
\begin{align*}
\mathcal L\mathbf Z_{I;b}
&=\mathcal L(D_a\mathbf Z_{I^-;b})\\
&=D_a(\mathcal L\mathbf Z_{I^-;b})
  +\mathcal C(\mathbf A_a,\mathbf Z_{I^-;b})\\
&=D_a\widetilde{\mathbf Q}_{I^-;b}
  +\mathcal C(\mathbf A_a,\mathbf Z_{I^-;b})\\
&=\widetilde{\mathbf Q}_{I;b}.
\end{align*}
Thus, induction on $|I|$ proves the identity
\[
 \mathcal L\mathbf{Z}_{I;b}=\widetilde{\mathbf{Q}}_{I;b}.
\]
This recurrence enumerates all summands. To make this explicit, consider formal expressions obtained by contracting, using \(\mathbf{g}\), factors \(\nabla^r\mathbf{R}\), \(\mathbf{T}\), \(\mathbf{Z}_{J;c}\), and \(D_t\mathbf{Z}_{J;c}\), leaving a single free vector index. Define a formal derivation \(\mathfrak d_a\) by the Leibniz rule on products and contractions and by the rules
\begin{align*}
 \mathfrak d_a\mathbf{T}&=D_t\mathbf{A}_a,
 &\mathfrak d_a\mathbf{Z}_{J;c}&=\mathbf{Z}_{(J,a);c},\\
 \mathfrak d_a(D_t\mathbf{Z}_{J;c})
 &=D_t\mathbf{Z}_{(J,a);c}+\mathbf{R}(\mathbf{A}_a,\mathbf{T})\mathbf{Z}_{J;c},
 &\mathfrak d_a(\nabla^r\mathbf{R})
 &=\nabla_{\mathbf{A}_a}(\nabla^r\mathbf{R})
  =(\nabla^{r+1}\mathbf{R})(\mathbf{A}_a,\,\cdot\,).
\end{align*}
In the last rule, $\cdot$ represents, in their original order, the $r$ differentiation arguments already present; the intrinsic indices of $\mathbf{R}$ remain within the tensor value. The new argument is placed first. The remaining tensor arguments are also differentiated, one at a time, by the Leibniz rule. Differentiating the contractions produces no terms because \(\nabla \mathbf{g}=0\). Let \(\mathfrak C_{a,J;b}\) be the six summands, with their coefficients, of
\begin{align*}
{}&-(\nabla_{\mathbf{T}}\mathbf{R})(\mathbf{A}_a,\mathbf{T})\mathbf{Z}_{J;b}
-\mathbf{R}(D_t\mathbf{A}_a,\mathbf{T})\mathbf{Z}_{J;b}
-2\mathbf{R}(\mathbf{A}_a,\mathbf{T})D_t\mathbf{Z}_{J;b}\\
&\qquad
-(\nabla_{\mathbf{A}_a}\mathbf{R})(\mathbf{Z}_{J;b},\mathbf{T})\mathbf{T}
-\mathbf{R}(\mathbf{Z}_{J;b},D_t\mathbf{A}_a)\mathbf{T}
-\mathbf{R}(\mathbf{Z}_{J;b},\mathbf{T})D_t\mathbf{A}_a .
\end{align*}
We now define finite multisets of expressions by
\[
 \mathfrak A_{\varnothing;b}:=\varnothing,
 \qquad
 \mathfrak A_{(J,a);b}
 :=\mathfrak d_a\mathfrak A_{J;b}\sqcup\mathfrak C_{a,J;b},
\]
where \(\mathfrak d_a\mathfrak A_{J;b}\) contains all summands produced by the preceding rules, retaining multiplicities and coefficients. Let us verify the resulting equality. For $I=\varnothing$, both sides are zero because $\widetilde{\mathbf Q}_{\varnothing;b}=0$ and $\mathfrak A_{\varnothing;b}=\varnothing$. Suppose that
\(\displaystyle \widetilde{\mathbf Q}_{J;b}
=\displaystyle\sum_{\mathfrak a\in\mathfrak A_{J;b}}\mathfrak a\)
and let $I=(J,a)$. Covariant differentiation $D_a$ acts on each product and contraction exactly according to the rules defining $\mathfrak d_a$. Hence,
\begin{align*}
\widetilde{\mathbf Q}_{I;b}
&=D_a\widetilde{\mathbf Q}_{J;b}
  +\mathcal C(\mathbf A_a,\mathbf Z_{J;b})=\displaystyle\sum_{\mathfrak a\in\mathfrak A_{J;b}}
  \mathfrak d_a\mathfrak a
  +\displaystyle\sum_{\mathfrak c\in\mathfrak C_{a,J;b}}\mathfrak c=\displaystyle\sum_{\mathfrak a\in\mathfrak A_{(J,a);b}}\mathfrak a.
\end{align*}
This proves by induction on the length of $I$ that
\begin{equation}
 \widetilde{\mathbf{Q}}_{I;b}=\displaystyle\sum_{\mathfrak a\in\mathfrak A_{I;b}}
 \mathfrak a.
 \label{eq:indexacion-exacta-fuentes-jacobi}
\end{equation}
Thus, \(\mathfrak A_{(J,a);b}\) is determined by \(\mathfrak A_{J;b}\) through the preceding union, which preserves the order of the letters and indexes every summand.

The indexing also has an exact grading. Assign weight zero to \(\mathbf{T}\) and weight \(\lvert J\rvert+1\) to every occurrence of \(\mathbf{Z}_{J;c}\) or \(D_t\mathbf{Z}_{J;c}\); in particular, \(\mathbf{A}_c=\mathbf{Z}_{\varnothing;c}\) has weight one. If \(\lvert I\rvert=k\geq1\), every summand of \eqref{eq:indexacion-exacta-fuentes-jacobi} satisfies
\begin{equation}
 \displaystyle\sum_{\text{occurrences of }\mathbf{Z}_{J;c}\text{ or }D_t\mathbf{Z}_{J;c}}
 (\lvert J\rvert+1)=k+1.
 \label{eq:peso-exacto-fuentes-jacobi}
\end{equation}
It contains one or more factors \(\nabla^r\mathbf{R}\), with \(0\leq r\leq k\), and each field factor contains at most one derivative \(D_t\). Moreover, every field appearing in the source has weight at most \(k\), strictly less than the weight \(k+1\) of the unknown \(\mathbf{Z}_{I;b}\). We verify the four assertions simultaneously. If $|I|=1$, write $I=(a)$. The six terms of $\mathfrak C_{a,\varnothing;b}$ contain the fields $\mathbf A_a$ and $\mathbf A_b$, or one of their derivatives $D_t$; their total weight is $1+1=2=|I|+1$. They contain a factor $\mathbf R$ or $\nabla\mathbf R$, contain no second time derivatives, and the fields appearing in them have weight one.

Now assume the assertions hold for a word $J$ of length $k-1$, and let $I=(J,a)$. In the terms arising from $\mathfrak d_a\mathfrak A_{J;b}$, the Leibniz rule makes $\mathfrak d_a$ act on one factor at a time. If it acts on $\mathbf Z_{L;c}$ or $D_t\mathbf Z_{L;c}$, it replaces that factor by one whose weight is exactly one greater; the additional curvature term in the rule for $D_t\mathbf Z_{L;c}$ also introduces $\mathbf A_a$, again of weight one. If it acts on $\mathbf T$, it introduces $D_t\mathbf A_a$, of weight one; and if it acts on $\nabla^r\mathbf R$, it raises the order of the curvature derivative from $r$ to $r+1$ and introduces the argument $\mathbf A_a$. In every case, the total weight increases by exactly one, from $k$ to $k+1$, no $D_t^2\mathbf Z$ appears, and the curvature order remains at most $k$. Finally, the six new terms of $\mathfrak C_{a,J;b}$ have weight $1+(|J|+1)=k+1$, curvature order at most one, and only first time derivatives. Since the field of greatest weight among them is $\mathbf Z_{J;b}$, of weight $k$, the entire source contains only fields of weight less than $k+1$. This completes the induction.

Restricting to \(\boldsymbol s=0\), set
\[
 \mathbf{Y}_{I;b}:=\mathbf{Z}_{I;b}|_{\boldsymbol s=0},
 \qquad
 \mathbf{Q}_{I;b}:=\widetilde{\mathbf{Q}}_{I;b}|_{\boldsymbol s=0}
 \in\Gamma(\gamma_v^*TM).
\]
Therefore,
\[
 \frac{D^2\mathbf{Y}_{I;b}}{dt^2}
 +\mathbf{R}(\mathbf{Y}_{I;b},\dot\gamma_v)\dot\gamma_v=\mathbf{Q}_{I;b}.
\]
All initial data are explicit. Since \(\alpha(\boldsymbol s,0)=p\),
\[
 \mathbf{A}_b(\boldsymbol s,0)=0,
 \qquad
 D_t\mathbf{A}_b(\boldsymbol s,0)=B_pw_b;
\]
and this, together with the same commutators, implies
\[
 \begin{array}{lll}
 \mathbf{Y}_{\varnothing;b}(0)=0,
 &\displaystyle \frac{D\mathbf{Y}_{\varnothing;b}}{dt}(0)=B_pw_b,&\\[6pt]
 \mathbf{Y}_{I;b}(0)=0,
 &\displaystyle \frac{D\mathbf{Y}_{I;b}}{dt}(0)=0,
 & |I|\geq1.
 \end{array}
\]
Indeed, at \(t=0\), all purely parametric derivatives of \(\mathbf{A}_b(\boldsymbol s,0)=0\) vanish, and every term arising when \(D_t\) is commuted with a nonempty word contains one of these fields; moreover, \(B_pw_b\) does not depend on \(\boldsymbol s\).

Restricting \eqref{eq:indexacion-exacta-fuentes-jacobi} preserves the grading \eqref{eq:peso-exacto-fuentes-jacobi}: it replaces \(\mathbf{T}\) by \(\dot\gamma_v\), \(\mathbf{Z}_{J;c}\) by \(\mathbf{Y}_{J;c}\), and \(D_t\mathbf{Z}_{J;c}\) by \(\frac{D\mathbf{Y}_{J;c}}{dt}\). Since \(|\dot\gamma_v|_{\mathbf{g}}=|B_pv|_{\mathbf{g}_p}
=\|v\|_{\mathbb R^n}\leq\rho\), the hypothesis on \(\nabla^r\mathbf{R}\) allows the bound to be expressed using $L^\infty([0,1])$ norms over the interval. For $k\in\mathbb N$, let $\displaystyle \mathcal W_{m,k}:=\bigcup_{\nu=0}^{k-1}\{0,\ldots,m\}^{\nu}$ be the set of words of length less than $k$, including the empty word. Define
\[
 \mathcal M_k
 :=
 \displaystyle\sum_{J\in\mathcal W_{m,k}}\displaystyle\sum_{c=0}^m
 \left(
 \|\mathbf Y_{J;c}\|_{L^\infty([0,1])}
 +
 \left\|\frac{D\mathbf Y_{J;c}}{dt}\right\|_{L^\infty([0,1])}
 \right),
\]
with the pointwise norms induced by $\mathbf g$. The index set of the sum is finite. By the grading already proved, each field factor in $\mathbf Q_{I;b}$ has weight at most $k$ when $|I|=k$. Therefore, there exist $C_k>0$ and $N_k\in\mathbb N$, independent of $p$, $v$, and the chosen coordinate directions, such that
\[
 \|\mathbf Q_{I;b}\|_{L^\infty([0,1])}
 \leq C_k(1+\mathcal M_k)^{N_k},
 \qquad |I|=k.
\]
Here $C_k$ and $N_k$ depend only on $n$, $k$, $\rho$, the bounds $\|\nabla^r\operatorname{Rm}_{\mathbf g}\|_\infty$ for $0\leq r\leq k$, and the finite number of contractions produced by the recurrence $\widetilde{\mathbf{Q}}_{(I^-,a);b}
=D_a\widetilde{\mathbf{Q}}_{I^-;b}
+\mathcal C(\mathbf{A}_a,\mathbf{Z}_{I^-;b})$; in particular, they do not depend on the center $p$.

To establish the bounds, fix \(v\) and use along \(\gamma_v\) the orthonormal frame obtained by parallel transport of \((B_pe_r)_{r=1}^n\). In this frame, the last equation has the form
\[
 y_{I;b}''(t)+\mathcal R_{p,v}(t)y_{I;b}(t)=q_{I;b}(t),
 \qquad
 \|\mathcal R_{p,v}(t)\|_{\operatorname{op},\mathbb R^n}\leq C_0\rho^2.
\]
It is equivalent to the first-order system
\[
 \frac d{dt}
 \begin{pmatrix}y_{I;b}\\y_{I;b}'\end{pmatrix}
 =
 \begin{pmatrix}0&I\\-\mathcal R_{p,v}(t)&0\end{pmatrix}
 \begin{pmatrix}y_{I;b}\\y_{I;b}'\end{pmatrix}
 +\begin{pmatrix}0\\q_{I;b}\end{pmatrix}.
\]
We proceed by induction on the weight. For weight one, $I=\varnothing$ and the equation is the homogeneous Jacobi equation with initial data $y_{\varnothing;b}(0)=0$ and $y_{\varnothing;b}'(0)=w_b$; the estimate obtained at the beginning of the proof provides a constant $B_1$ such that $\mathcal M_1\leq B_1$.

Suppose that, for some $k\geq1$, all fields of weight at most $k$ satisfy the corresponding bounds, and hence $\mathcal M_k\leq B_k$. If $|I|=k$, the preceding estimate gives
\[
 \|q_{I;b}\|_{L^\infty([0,1])}
 \leq C_k(1+B_k)^{N_k}=:Q_k.
\]
Let $z_{I;b}=(y_{I;b},y_{I;b}')$. The matrix of the first-order system has norm at most $L:=1+C_0\rho^2$. Since $z_{I;b}(0)=0$ for $k\geq1$, the integral inequality associated with the system is
\[
 |z_{I;b}(t)|
 \leq
 L\int_0^t|z_{I;b}(\tau)|\,d\tau+tQ_k.
\]
Grönwall's Lemma~\ref{lema: gronwall} implies
\[
 \|y_{I;b}\|_{L^\infty([0,1])}
 +\|y_{I;b}'\|_{L^\infty([0,1])}
 \leq 2Q_ke^L.
\]
There are only finitely many words $I$ and letters $b$, their number depending on $k$ and $m$. Adding these bounds to those forming $\mathcal M_k$ gives a constant $B_{k+1}$ with $\mathcal M_{k+1}\leq B_{k+1}$. This constant depends only on the geometric bounds listed above and on $B_1,\dots,B_k$, never on $p$, $v$, or the chosen letters. This completes the induction on the weight.

It remains to pass from these covariant derivatives to ordinary derivatives of the matrices. To this end, over the full variation, let \((\mathbf{E}_r(\boldsymbol s,t))_{r=1}^n\) be the frame determined by
\[
 D_t\mathbf{E}_r=0,
 \qquad
 \mathbf{E}_r(\boldsymbol s,0)=B_pe_r.
\]
It is smooth in \(\boldsymbol s\) and orthonormal for every \(\boldsymbol s\). If \(\mathbf{W}_{I;r}:=D_I\mathbf{E}_r\), then
\[
 D_t\mathbf{W}_{\varnothing;r}=0,
 \qquad
 \mathbf{W}_{\varnothing;r}(0)=B_pe_r,
\]
and, defining \(\mathbf{S}_{\varnothing;r}:=0\) and, for \(I=(I^-,a)\),
\[
 \mathbf{S}_{I;r}:=D_a\mathbf{S}_{I^-;r}-\mathbf{R}(\mathbf{A}_a,\mathbf{T})\mathbf{W}_{I^-;r},
\]
the commutator fixed above gives precisely
\[
 D_t\mathbf{W}_{I;r}=\mathbf{S}_{I;r},
 \qquad
 \mathbf{W}_{I;r}(0)=0\quad (|I|\geq1).
\]
After expanding \(D_a\mathbf{S}_{I^-;r}\), its right-hand side contains only curvature, the variational fields already bounded, and \(\mathbf{W}_{J;r}\) with \(|J|<|I|\). For $I=\varnothing$, $\mathbf W_{\varnothing;r}=\mathbf E_r$ is parallel and unit length, so $|\mathbf W_{\varnothing;r}|_{\mathbf g}=1$. If $I=(a)$, then
\[
 \mathbf S_{(a);r}
 =-\mathbf R(\mathbf A_a,\mathbf T)\mathbf E_r,
 \qquad
 \mathbf W_{(a);r}(t)=\int_0^t\mathbf S_{(a);r}(\tau)\,d\tau,
\]
where the integral is taken after identifying the fibers by parallel transport; the bounds on $\mathbf R$, $\mathbf A_a$, and $\mathbf T$ bound both sides uniformly. Now suppose all $\mathbf W_{J;r}$ with $|J|<k$ are bounded, and let $|I|=k$. The full Leibniz expansion of $\mathbf S_{I;r}$ contains only these lower-order fields, curvature derivatives of order at most $k-1$, and the variational fields already bounded. Therefore, there exists $C_k$ such that $\|\mathbf S_{I;r}\|_{L^\infty([0,1])}\leq C_k$. Since $\mathbf W_{I;r}(0)=0$,
\[
 \mathbf W_{I;r}(t)=\int_0^t\mathbf S_{I;r}(\tau)\,d\tau,
 \qquad
 \|\mathbf W_{I;r}\|_{L^\infty([0,1])}\leq C_k.
\]
This completes the induction.

Now take \(w_0=e_j\). If
\[
 y_j^r(v,t):=\langle \mathbf{J}_j(t),\mathbf{E}_r(v,t)\rangle_{\mathbf{g}},
\]
the variations with \(w_1,\ldots,w_k\) compute its ordinary derivatives in \(v\). More precisely, for an ordered word \(I=(a_1,\ldots,a_k)\) and a subset \(S\subseteq\{1,\ldots,k\}\), denote by \(I_S\) the subword retaining the original order. Metric compatibility of the connection and the Leibniz rule give the exact triangular identity
\[
 \partial_I y_j^r(v,t)
 =\displaystyle\sum_{S\subseteq\{1,\ldots,k\}}
 \left\langle
   \mathbf{Y}_{I_S;0}(t),\mathbf{W}_{I_{S^c};r}(t)
 \right\rangle_{\mathbf{g}},
\]
where \(\partial_I=\partial_{s_{a_k}}\cdots\partial_{s_{a_1}}|_{\boldsymbol s=0}\). The term corresponding to \(S=\{1,\ldots,k\}\) is \(\langle \mathbf{Y}_{I;0},\mathbf{E}_r\rangle_{\mathbf{g}}\), and all the other terms contain fewer derivatives on \(\mathbf{A}_0\). The bounds obtained for \(\mathbf{Y}_{I;0}\) and \(\mathbf{W}_{I;r}\) therefore imply uniform bounds for all ordinary derivatives of \(y_j^r(v,t)\).

At \(t=1\), the matrix \(\mathbf{Y}(v)=(y_j^r(v,1))_{rj}\) represents \(d\Psi_p(v)\) in the parallel frame, and
\[
 g_{ij}^{(p)}(v)
 =\langle \mathbf{J}_i(1),\mathbf{J}_j(1)\rangle_{\mathbf{g}}
 =\displaystyle\sum_{r=1}^n y_i^r(v,1)y_j^r(v,1).
\]
The scalar Leibniz rule thus gives, for every multi-index \(\alpha\),
\[
 \displaystyle\sup_{p\in M}\displaystyle\sup_{v\in B(0,\rho)}
 |D^\alpha g_{ij}^{(p)}(v)|<\infty.
\]
The uniform ellipticity already established gives a uniform bound for \(g_{(p)}^{ij}\). For \(|\alpha|\geq1\), differentiate the identity \(\displaystyle \sum_{\ell=1}^n g_{(p)}^{i\ell}g_{\ell j}^{(p)}=\delta_j^i\), isolate the term containing \(\beta=\alpha\), and multiply by the inverse metric. Componentwise, this gives
\[
 D^\alpha g_{(p)}^{ik}
 =-\displaystyle\sum_{\substack{\beta\in\mathbb N_0^n\\\beta\leq\alpha,\ \beta\neq\alpha}}
 \binom{\alpha}{\beta}
 \displaystyle\sum_{\ell,j=1}^n
 (D^\beta g_{(p)}^{i\ell})
 (D^{\alpha-\beta}g_{\ell j}^{(p)})g_{(p)}^{jk}.
\]
Induction on $|\alpha|$ is now direct. For $|\alpha|=0$, we already have the uniform bound for $g_{(p)}^{ij}$. Suppose all derivatives $D^\beta g_{(p)}^{ij}$ with $|\beta|<k$ are bounded, and let $|\alpha|=k$. In every summand of the preceding formula, $|\beta|<k$ and $1\leq|\alpha-\beta|\leq k$; by the induction hypothesis, the bounds already proved for $D^{\alpha-\beta}g_{\ell j}^{(p)}$, and the bound for $g_{(p)}^{jk}$, every product is uniformly bounded. The sum contains finitely many terms, their number depending only on $\alpha$. Thus, $D^\alpha g_{(p)}^{ik}$ is also uniformly bounded, and induction proves the assertion for all derivatives of the components $g_{(p)}^{ij}$.

We have proved (2) with \(r_0:=\rho\). For \(0<r\leq r_0\), the bounds in charts of radius \(r\) follow by restricting the charts of radius \(r_0\).

We now prove the additional assertion about coordinate transitions. Fix \(0<r\leq r_0\). On an overlap, set \(\Phi=\phi_j^{\operatorname{geo}}\circ
(\phi_i^{\operatorname{geo}})^{-1}\). The tensor identity
\[
g_{ab}^{(i)}
=\displaystyle\sum_{\ell,m=1}^n
 (g_{\ell m}^{(j)}\circ\Phi)
 (\partial_a\Phi^\ell)(\partial_b\Phi^m)
\]
and uniform ellipticity bound \(D\Phi\): indeed, evaluating the equality on a vector \(\xi\), the lower bound for \(g^{(j)}\) and the upper bound for \(g^{(i)}\) give \(c|D\Phi\,\xi|^2\leq C|\xi|^2\). Interchanging \(i\) and \(j\) bounds the differential of the inverse transition. The Christoffel symbols and all their derivatives are uniformly bounded because their components are universal expressions in \(g_{(i)}^{ab}\) and the derivatives of \(g_{ab}^{(i)}\). The transformation law for the connection reads
\[
\partial_a\partial_b\Phi^k
=\Gamma_{i,ab}^{c}\,\partial_c\Phi^k
-\Gamma_{j,\ell m}^{k}(\Phi)
\partial_a\Phi^\ell\partial_b\Phi^m.
\]
This equality bounds \(D^2\Phi\). For the inductive step, suppose that, for some \(m\geq2\), all derivatives of \(\Phi\) of order at most \(m\) are uniformly bounded, and let \(\alpha\in\mathbb N_0^n\) with \(\lvert\alpha\rvert=m-1\). Apply \(D^\alpha\) to the preceding identity. The Leibniz rule expresses the derivative of the first summand in terms of derivatives of \(\Gamma_i\) of order at most \(m-1\) and derivatives of \(\Phi\) of order at most \(m\).

For the second summand, first apply the Leibniz rule to its three factors. In the composite factor, set \(A:=\Gamma_{j,\ell m}^{k}\). The Faà di Bruno formula of Theorem~\ref{faa di bruno multivariable}, written in the same partial-derivative notation used here, states that, for \(1\leq\lvert\mu\rvert\leq m-1\),
\begin{equation}
\label{eq:faa-di-bruno-geometria-acotada}
D^\mu(A\circ\Phi)(x)
=
\displaystyle\sum_{\substack{
\beta\in\mathbb N_0^n,\ 1\leq\lvert\beta\rvert\leq\lvert\mu\rvert\\
\gamma=(\gamma_1,\dots,\gamma_{\lvert\beta\rvert}),\quad
\lvert\gamma_\nu\rvert>0\\
\gamma_1+\cdots+\gamma_{\lvert\beta\rvert}=\mu\\
l=(l_1,\dots,l_{\lvert\beta\rvert}),\quad
l_\nu\in\{1,\dots,n\}\\
\#\{\nu\mid l_\nu=a\}=\beta_a,\quad 1\leq a\leq n}}
c_{\mu,\beta,\gamma,l}\,
(D^\beta A)(\Phi(x))
\prod_{\nu=1}^{\lvert\beta\rvert}
D^{\gamma_\nu}\Phi^{l_\nu}(x).
\end{equation}
In each term, \(\lvert\beta\rvert\leq\lvert\mu\rvert\leq m-1\), \(\#\{\nu\mid l_\nu=a\}=\beta_a\) for \(1\leq a\leq n\), and \(1\leq\lvert\gamma_\nu\rvert\leq\lvert\mu\rvert\leq m-1\). Thus, the derivatives of \(A\) are bounded by the metric bounds, and the derivatives of \(\Phi\) appearing in \eqref{eq:faa-di-bruno-geometria-acotada} are bounded by the induction hypothesis. The other two factors are first derivatives of \(\Phi\); after applying the Leibniz rule, they involve derivatives of order at most \(m\). Every term on the right-hand side is therefore bounded by a constant independent of the charts. The left-hand side is \(D^\alpha\partial_a\partial_b\Phi^k\), of order \(m+1\), which completes the induction: every multi-index of length $m+1\geq2$ can be written as $\alpha+e_a+e_b$, allowing $a=b$. Interchanging \(i\) and \(j\) gives the same estimate for the inverse transition.

It remains to prove that (2) implies (1). Fix \(x\in M\) and use the normal chart of radius \(r_0\) centered at \(x\). The formula
\[
\Gamma^\ell_{ij}
=
\frac12\displaystyle\sum_{a=1}^n g^{\ell a}
\bigl(\partial_i g_{ja}+\partial_j g_{ia}-\partial_a g_{ij}\bigr)
\]
gives, for each multi-index \(\alpha\), the exact identity
\[
\partial^\alpha\Gamma^\ell_{ij}
=
\frac12\displaystyle\sum_{a=1}^n\displaystyle\sum_{\substack{\beta\in\mathbb N_0^n\\\beta\leq\alpha}}\binom{\alpha}{\beta}
(\partial^\beta g^{\ell a})
\left(
\partial^{\alpha-\beta+e_i}g_{ja}
+\partial^{\alpha-\beta+e_j}g_{ia}
-\partial^{\alpha-\beta+e_a}g_{ij}
\right).
\]
The sum over \(a\) and \(\beta\) is finite, and condition (2) bounds each of its factors uniformly.

The curvature components satisfy
\[
R^\ell{}_{ijk}
=\partial_i\Gamma^\ell_{jk}-\partial_j\Gamma^\ell_{ik}
+\displaystyle\sum_{m=1}^n\left(\Gamma^m_{jk}\Gamma^\ell_{im}
-\Gamma^m_{ik}\Gamma^\ell_{jm}\right).
\]
This formula and the Leibniz rule give uniform bounds for all partial derivatives of \(R\). Finally, the coordinate formula for the covariant derivative of a tensor consists of its partial derivative plus a finite sum of contractions with Christoffel symbols. We proceed by induction on $q$. For $q=0$, the components of $R$ and all their partial derivatives are uniformly bounded by the preceding calculation. Suppose that each component of $\nabla^qR$ is a finite universal sum of products of partial derivatives of $R$ and partial derivatives of the Christoffel symbols, all of order at most $q$. Write the curvature tensor as a tensor of type $(0,4)$, lowering its contravariant index using $\mathbf g$. Lemma~\ref{lem:local-expression-higher-order} gives, for $I=(i_1,\ldots,i_{q+4})\in\{1,\ldots,n\}^{q+4}$ and $a\in\{1,\ldots,n\}$,
\[
 (\nabla_a\nabla^qR)_{i_1\cdots i_{q+4}}
 =\partial_a(\nabla^qR)_{i_1\cdots i_{q+4}}
 -\displaystyle\sum_{\nu=1}^{q+4}\displaystyle\sum_{c=1}^n
   \Gamma^c_{a i_\nu}
   (\nabla^qR)_{i_1\cdots i_{\nu-1}c i_{\nu+1}\cdots i_{q+4}}.
\]
When differentiating the products in the induction hypothesis, the Leibniz rule increases by one the order of exactly one factor; the connection terms add a factor $\Gamma$ without increasing the other orders. Thus, each component of $\nabla^{q+1}R$ has the same form with orders at most $q+1$. All these components are uniformly bounded, completing the induction. At the center of the normal chart, \(g_{ij}(x)=\delta_{ij}\); consequently, the tensor norm of \(\nabla^q\mathbf{R}(x)\) is bounded by a constant independent of \(x\). Since \(x\) was arbitrary, this gives (1) and concludes the equivalence.
\end{proof}
To justify the calculation of the injectivity radius for the nonpositive-curvature models appearing below, we record the required global result.

\begin{theorem}[Cartan--Hadamard]\label{teo:cartan-hadamard}
Let $(M,\mathbf{g})$ be a complete, connected Riemannian manifold without boundary whose sectional curvature is nonpositive. For every $p\in M$, the exponential map
\[
\exp_p\colon T_pM\longrightarrow M
\]
is a covering map. If, in addition, $M$ is simply connected, then $\exp_p$ is a global diffeomorphism and $\operatorname{inj}(p)=+\infty$.
\end{theorem}

\begin{proof}
Let $\gamma$ be a geodesic starting at $p$, and let $\mathbf{J}$ be a Jacobi field perpendicular to $\dot\gamma$, with $\mathbf{J}(0)=0$. The Jacobi equation and the curvature hypothesis give
\[
\frac{1}{2}\frac{d^2}{dt^2}\|\mathbf{J}(t)\|_{\mathbf{g}}^2
=\left\|\frac{D\mathbf{J}}{dt}(t)\right\|_{\mathbf{g}}^2
-\bigl\langle \mathbf{R}(\mathbf{J},\dot\gamma)\dot\gamma,\mathbf{J}\bigr\rangle_{\mathbf{g}}
\geq 0.
\]
If $\mathbf{J}$ vanished again, the nonnegative convex function $\|\mathbf{J}\|_{\mathbf{g}}^2$ would vanish between the two endpoints; uniqueness for the Jacobi equation would imply $\mathbf{J}\equiv0$. Thus, there are no conjugate points and $d\exp_p(v)$ is invertible for every $v\in T_pM$; hence, $\exp_p$ is a local diffeomorphism.

To make the differential estimate precise, decompose the initial datum of a Jacobi field with \(\mathbf{J}(0)=0\) into its components parallel and perpendicular to \(\dot\gamma(0)\). If \(P_t\) denotes parallel transport along \(\gamma\), the parallel component is exactly \(\mathbf{J}_\parallel(t)
=tP_t(D_t\mathbf{J}(0)_\parallel)\). For the perpendicular component, at points where \(f(t):=\lVert \mathbf{J}(t)\rVert_{\mathbf{g}}\) does not vanish, the Jacobi equation gives
\[
 f''(t)
 =
 \frac{
 \lVert D_t\mathbf{J}\rVert_{\mathbf{g}}^2f(t)^2
 -\langle D_t\mathbf{J},\mathbf{J}\rangle_{\mathbf{g}}^2
 -\langle \mathbf{R}(\mathbf{J},\dot\gamma)\dot\gamma,\mathbf{J}\rangle_{\mathbf{g}} f(t)^2
 }{f(t)^3}
 \geq0.
\]
Here we used Cauchy--Schwarz and the nonpositive sectional curvature. Since \(f(0)=0\) and \(f'(0^+)=\lVert D_t\mathbf{J}(0)\rVert_{\mathbf{g}}\), convexity implies \(f(t)\geq t\lVert D_t\mathbf{J}(0)\rVert_{\mathbf{g}}\). The parallel and perpendicular components remain orthogonal. Evaluating at \(t=1\) gives
\[
\|d\exp_p(v)w\|_{\mathbf{g}}\geq \|w\|_{\mathbf{g}_p},
\qquad v,w\in T_pM.
\]
Consequently, the metric $\exp_p^*\mathbf{g}$ dominates the Euclidean metric of $T_pM$ and is complete: a Cauchy sequence for $\exp_p^*\mathbf{g}$ is Cauchy in $T_pM$. If $v$ is its limit, continuity and positivity of $\exp_p^*\mathbf{g}$ provide a neighborhood $V$ of $v$ and constants $c_v,C_v>0$ such that
\[
 c_v\|w\|_{\mathbf{g}_p}^2
 \leq (\exp_p^*\mathbf{g})_z(w,w)
 \leq C_v\|w\|_{\mathbf{g}_p}^2,
 \qquad z\in V,\quad w\in T_pM.
\]
The sequence then converges for the distance of $\exp_p^*\mathbf{g}$ as well.

The exponential map is surjective by Hopf--Rinow (Theorem~\ref{teo:variedades-riemannianas-hopf-rinow}). We record the passage from local isometry to covering map. Let \(y\in M\), and choose \(\varepsilon>0\) so that \(\exp_y\) is a diffeomorphism on the ball of radius \(\varepsilon\). For each \(x\in\exp_p^{-1}(y)\), completeness of \((T_pM,\exp_p^*\mathbf{g})\) allows us to define
\[
s_x(\exp_y v)
:=
\exp_x^{\,\exp_p^*\mathbf{g}}\bigl((d\exp_p)_x^{-1}v\bigr),
\qquad |v|_{\mathbf{g}_y}<\varepsilon.
\]
A local isometry maps geodesics to geodesics and preserves their initial data; geodesic uniqueness proves \(\exp_p\circ s_x=\operatorname{id}_{B_{\mathbf{g}}(y,\varepsilon)}\). The images of \(s_x\) are disjoint: if two contained the same point, the lifts of the reversed radial geodesic, with the same initial datum determined by \(d\exp_p\), would coincide and have the same endpoint over \(y\). Finally, if \(q\in\exp_p^{-1}(B_{\mathbf{g}}(y,\varepsilon))\), lift from \(q\) the radial geodesic ending at \(y\); completeness extends the lift to its endpoint \(x\in\exp_p^{-1}(y)\), and the reversed geodesic shows that \(q=s_x(\exp_p(q))\). Therefore,
\[
\exp_p^{-1}(B_{\mathbf{g}}(y,\varepsilon))
=\coprod_{x\in\exp_p^{-1}(y)}
s_x(B_{\mathbf{g}}(y,\varepsilon)),
\]
and each summand projects diffeomorphically onto the ball. This proves that \(\exp_p\) is a covering map. Since the finite-dimensional vector space $T_pM$ is simply connected, it is the universal covering space; if $M$ is also simply connected, the covering has only one sheet. The last assertion follows from Definition~\ref{def: radio de inyectividad}.
\end{proof}

\begin{example}\label{ej:geometria-acotada-presentamos-algunos-ejemplos-fundamentales-variedades-riemannianas}
We present some fundamental examples of smooth Riemannian manifolds without boundary and with bounded geometry, explicitly verifying the conditions of Definition~\ref{def: geometria acotada}.

\begin{enumerate}
\item Euclidean space $(\mathbb{R}^{n},\overline{\mathbf{g}})$ has bounded geometry.

Work in the standard coordinates of $\mathbb{R}^{n}$. In these coordinates,
$
\overline g_{ij}=\delta_{ij},
$
so
$
\partial_{k}\overline g_{ij}=0
$
for all $i,j,k\in\{1,\dots,n\}$.

By Corollary~\ref{expresion local christoffel levi civita}, the Christoffel symbols of the Levi--Civita connection satisfy
$
\Gamma_{ij}^{k}=0
$
for all $i,j,k\in\{1,\dots,n\}$.

Substituting into Proposition~\ref{1,3 curvatura en coordenadas} gives
$
R_{ijk}{}^{l}=0
$
for all $i,j,k,l\in\{1,\dots,n\}$, and consequently
$
\mathbf{R}=0.
$
Therefore,
$
\nabla^{m}\mathbf{R}=0,\; \forall m\in\mathbb{N}_{0}.
$

We now study the exponential map. Let $p\in\mathbb{R}^{n}$ and $v\in T_{p}\mathbb{R}^{n}\cong\mathbb{R}^{n}$. In standard coordinates, the geodesic equation is
$
\displaystyle\frac{d^{2}x^{k}}{dt^{2}}(t)=0,
\; k\in\{1,\dots,n\},
$
since $\Gamma_{ij}^{k}=0$.

Integrating twice shows that every solution has the form
$
x(t)=a+tb,
$
with $a,b\in\mathbb{R}^{n}$. Imposing the initial conditions
$
x(0)=p,
\;
\displaystyle\frac{dx}{dt}(0)=v,
$
gives $a=p$ and $b=v$, so the geodesic with initial data $(p,v)$ is given by
$
\gamma_{v}(t)=p+tv.
$

By definition of the exponential map,
$
\exp_{p}(v)=\gamma_{v}(1),
$
and therefore
$
\exp_{p}(v)=p+v.
$

In particular, the exponential map
$
\exp_{p}\colon \mathbb{R}^{n}\longrightarrow \mathbb{R}^{n}
$
is a global diffeomorphism, with inverse given by $q\mapsto q-p$. Consequently,
$
\operatorname{inj}(p)=+\infty
$
for every $p\in\mathbb{R}^{n}$, and therefore
$
\operatorname{inj}(\mathbb{R}^{n})=+\infty.
$

We conclude that $(\mathbb{R}^{n},\overline{\mathbf{g}})$ has bounded geometry.

\item Consider the upper half-plane
\[
\mathbb{H}^{2}:=\{(x,y)\in \mathbb{R}^{2}\mid y>0\}
\]
equipped with the Riemannian metric
\[
\mathbf{g}=\frac{(\mathbf{d}x)^{2}+(\mathbf{d}y)^{2}}{y^{2}}.
\]
Work in the global coordinates induced by the inclusion $\mathbb{H}^{2}\subseteq \mathbb{R}^{2}$, namely $(x^{1},x^{2})=(x,y)$. Denote the coordinate fields by $\boldsymbol{\partial}_i:=\partial/\partial x^i$.

The coefficients of the metric and its inverse are given by
$
g_{ij}=\displaystyle\frac{1}{y^{2}}\delta_{ij}
$
and
$
g^{ij}=y^{2}\delta_{ij}.
$
In particular,
$
g_{11}=g_{22}=\displaystyle\frac{1}{y^{2}}
$
and
$
g_{12}=g_{21}=0.
$

Since $g_{ij}$ depends only on the variable $y=x^{2}$,
$
\partial_{1}g_{ij}=0
$
for all $i,j\in\{1,2\}$, whereas
$
\partial_{2}g_{11}=\partial_{2}g_{22}=-\displaystyle\frac{2}{y^{3}}
$
and
$
\partial_{2}g_{12}=\partial_{2}g_{21}=0.
$

Using Corollary~\ref{expresion local christoffel levi civita}, we obtain
\[
\Gamma^{k}_{ij}
=
\frac{1}{2}g^{k\ell}
\bigl(
\partial_{i}g_{j\ell}
+
\partial_{j}g_{i\ell}
-
\partial_{\ell}g_{ij}
\bigr).
\]

First compute the symbols with upper index $1$. Since $g^{11}=y^{2}$ and $g^{12}=0$, we have
\[
\Gamma^{1}_{11}
=
\frac{1}{2}g^{1\ell}
\bigl(
\partial_{1}g_{1\ell}
+
\partial_{1}g_{1\ell}
-
\partial_{\ell}g_{11}
\bigr)
=0.
\]
Moreover,
\[
\Gamma^{1}_{12}
=
\frac{1}{2}g^{1\ell}
\bigl(
\partial_{1}g_{2\ell}
+
\partial_{2}g_{1\ell}
-
\partial_{\ell}g_{12}
\bigr)
=
\frac{1}{2}g^{11}\partial_{2}g_{11}
=
\frac{1}{2}y^{2}\left(-\frac{2}{y^{3}}\right)
=
-\frac{1}{y}.
\]
By symmetry of the Levi--Civita connection,
$
\Gamma^{1}_{21}=\Gamma^{1}_{12}=-\displaystyle\frac{1}{y}.
$
Finally,
\[
\Gamma^{1}_{22}
=
\frac{1}{2}g^{1\ell}
\bigl(
\partial_{2}g_{2\ell}
+
\partial_{2}g_{2\ell}
-
\partial_{\ell}g_{22}
\bigr)
=0.
\]

Now compute the symbols with upper index $2$. Since $g^{22}=y^{2}$ and $g^{21}=0$, we obtain
\[
\Gamma^{2}_{11}
=
\frac{1}{2}g^{2\ell}
\bigl(
\partial_{1}g_{1\ell}
+
\partial_{1}g_{1\ell}
-
\partial_{\ell}g_{11}
\bigr)
=
-\frac{1}{2}g^{22}\partial_{2}g_{11}
=
-\frac{1}{2}y^{2}\left(-\frac{2}{y^{3}}\right)
=
\frac{1}{y}.
\]
Similarly,
\[
\Gamma^{2}_{12}
=
\frac{1}{2}g^{2\ell}
\bigl(
\partial_{1}g_{2\ell}
+
\partial_{2}g_{1\ell}
-
\partial_{\ell}g_{12}
\bigr)
=0,
\]
and, by symmetry,
$
\Gamma^{2}_{21}=0.
$
Lastly,
\[
\Gamma^{2}_{22}
=
\frac{1}{2}g^{2\ell}
\bigl(
\partial_{2}g_{2\ell}
+
\partial_{2}g_{2\ell}
-
\partial_{\ell}g_{22}
\bigr)
=
\frac{1}{2}g^{22}\partial_{2}g_{22}
=
\frac{1}{2}y^{2}\left(-\frac{2}{y^{3}}\right)
=
-\frac{1}{y}.
\]

Thus,
\[
\Gamma^{1}_{11}=0,
\qquad
\Gamma^{1}_{12}=\Gamma^{1}_{21}=-\frac{1}{y},
\qquad
\Gamma^{1}_{22}=0,
\]
\[
\Gamma^{2}_{11}=\frac{1}{y},
\qquad
\Gamma^{2}_{12}=\Gamma^{2}_{21}=0,
\qquad
\Gamma^{2}_{22}=-\frac{1}{y}.
\]

We now compute the curvature tensor using Proposition~\ref{1,3 curvatura en coordenadas}:
\[
R_{ijk}{}^{\ell}
=
\partial_{i}\Gamma^{\ell}_{jk}
-
\partial_{j}\Gamma^{\ell}_{ik}
+
\Gamma^{m}_{jk}\Gamma^{\ell}_{im}
-
\Gamma^{m}_{ik}\Gamma^{\ell}_{jm}.
\]

Consider the component $R_{122}{}^{1}$. We have
\[
\partial_{1}\Gamma^{1}_{22}=0,
\qquad
\partial_{2}\Gamma^{1}_{12}
=
\partial_{2}\left(-\frac{1}{y}\right)
=
\frac{1}{y^{2}}.
\]
Moreover,
\[
\Gamma^{m}_{22}\Gamma^{1}_{1m}
=
\Gamma^{1}_{22}\Gamma^{1}_{11}
+
\Gamma^{2}_{22}\Gamma^{1}_{12}
=
0+\left(-\frac{1}{y}\right)\left(-\frac{1}{y}\right)
=
\frac{1}{y^{2}},
\]
and
\[
\Gamma^{m}_{12}\Gamma^{1}_{2m}
=
\Gamma^{1}_{12}\Gamma^{1}_{21}
+
\Gamma^{2}_{12}\Gamma^{1}_{22}
=
\left(-\frac{1}{y}\right)\left(-\frac{1}{y}\right)+0
=
\frac{1}{y^{2}}.
\]

Therefore,
\[
R_{122}{}^{1}
=
0-\frac{1}{y^{2}}+\frac{1}{y^{2}}-\frac{1}{y^{2}}
=
-\frac{1}{y^{2}}.
\]

Lowering indices gives
\[
R_{1221}
=
g_{1\ell}R_{122}{}^{\ell}
=
g_{11}R_{122}{}^{1}
=
\frac{1}{y^{2}}\left(-\frac{1}{y^{2}}\right)
=
-\frac{1}{y^{4}}.
\]

We now compute the sectional curvature of the plane spanned by $\boldsymbol{\partial}_{1}$ and $\boldsymbol{\partial}_{2}$. By definition,
\[
K(\boldsymbol{\partial}_{1},\boldsymbol{\partial}_{2})
=
\frac{\langle \mathbf{R}(\boldsymbol{\partial}_{1},\boldsymbol{\partial}_{2})
\boldsymbol{\partial}_{2},\boldsymbol{\partial}_{1}\rangle}
{\langle \boldsymbol{\partial}_{1},\boldsymbol{\partial}_{1}\rangle
\langle \boldsymbol{\partial}_{2},\boldsymbol{\partial}_{2}\rangle
-\langle \boldsymbol{\partial}_{1},\boldsymbol{\partial}_{2}\rangle^{2}}.
\]
The denominator is
\[
\langle \boldsymbol{\partial}_{1},\boldsymbol{\partial}_{1}\rangle
\langle \boldsymbol{\partial}_{2},\boldsymbol{\partial}_{2}\rangle
-\langle \boldsymbol{\partial}_{1},\boldsymbol{\partial}_{2}\rangle^{2}
=
g_{11}g_{22}-(g_{12})^{2}
=
\frac{1}{y^{4}},
\]
whereas the numerator is
\[
\langle \mathbf{R}(\boldsymbol{\partial}_{1},\boldsymbol{\partial}_{2})
\boldsymbol{\partial}_{2},\boldsymbol{\partial}_{1}\rangle
=
R_{1221}
=
-\frac{1}{y^{4}}.
\]
Consequently,
\[
K(\boldsymbol{\partial}_{1},\boldsymbol{\partial}_{2})
=
\frac{-\frac{1}{y^{4}}}{\frac{1}{y^{4}}}
=
-1.
\]

We have proved that the sectional curvature is constant and equal to $-1$. By Theorem~\ref{thm: tensor curvatura constante}, it follows that
\[
\mathbf{R}(\mathbf{X},\mathbf{Y})\mathbf{Z}
=
-\bigl(\mathbf{g}(\mathbf{Y},\mathbf{Z})\mathbf{X}-\mathbf{g}(\mathbf{X},\mathbf{Z})\mathbf{Y}\bigr).
\]

We now verify that $\nabla \mathbf{R}=0$. By definition of the covariant derivative of a tensor of type $(1,3)$,
\[
(\nabla_{\mathbf{W}}\mathbf{R})(\mathbf{X},\mathbf{Y})\mathbf{Z}
=
\nabla_{\mathbf{W}}\bigl(\mathbf{R}(\mathbf{X},\mathbf{Y})\mathbf{Z}\bigr)
-
\mathbf{R}(\nabla_{\mathbf{W}}\mathbf{X},\mathbf{Y})\mathbf{Z}
-
\mathbf{R}(\mathbf{X},\nabla_{\mathbf{W}}\mathbf{Y})\mathbf{Z}
-
\mathbf{R}(\mathbf{X},\mathbf{Y})\nabla_{\mathbf{W}}\mathbf{Z}.
\]

Substituting the preceding expression for $\mathbf{R}(\mathbf{X},\mathbf{Y})\mathbf{Z}$ gives
\[
\nabla_{\mathbf{W}}\bigl(\mathbf{R}(\mathbf{X},\mathbf{Y})\mathbf{Z}\bigr)
=
-\nabla_{\mathbf{W}}\bigl(\mathbf{g}(\mathbf{Y},\mathbf{Z})\mathbf{X}-\mathbf{g}(\mathbf{X},\mathbf{Z})\mathbf{Y}\bigr).
\]

Applying the Leibniz rule for the connection gives
\[
\nabla_{\mathbf{W}}\bigl(\mathbf{g}(\mathbf{Y},\mathbf{Z})\mathbf{X}\bigr)
=
\mathbf{W}\bigl(\mathbf{g}(\mathbf{Y},\mathbf{Z})\bigr)\mathbf{X} + \mathbf{g}(\mathbf{Y},\mathbf{Z})\nabla_{\mathbf{W}}\mathbf{X},
\]
\[
\nabla_{\mathbf{W}}\bigl(\mathbf{g}(\mathbf{X},\mathbf{Z})\mathbf{Y}\bigr)
=
\mathbf{W}\bigl(\mathbf{g}(\mathbf{X},\mathbf{Z})\bigr)\mathbf{Y} + \mathbf{g}(\mathbf{X},\mathbf{Z})\nabla_{\mathbf{W}}\mathbf{Y}.
\]
Therefore,
\[
\nabla_{\mathbf{W}}\bigl(\mathbf{R}(\mathbf{X},\mathbf{Y})\mathbf{Z}\bigr)
=
-\Bigl(
\mathbf{W}(\mathbf{g}(\mathbf{Y},\mathbf{Z}))\mathbf{X} + \mathbf{g}(\mathbf{Y},\mathbf{Z})\nabla_{\mathbf{W}}\mathbf{X}
-
\mathbf{W}(\mathbf{g}(\mathbf{X},\mathbf{Z}))\mathbf{Y} - \mathbf{g}(\mathbf{X},\mathbf{Z})\nabla_{\mathbf{W}}\mathbf{Y}
\Bigr).
\]

Now using compatibility of the Levi--Civita connection with the metric,
\[
\mathbf{W}(\mathbf{g}(\mathbf{Y},\mathbf{Z}))=\mathbf{g}(\nabla_{\mathbf{W}}\mathbf{Y},\mathbf{Z})+\mathbf{g}(\mathbf{Y},\nabla_{\mathbf{W}}\mathbf{Z}),
\qquad
\mathbf{W}(\mathbf{g}(\mathbf{X},\mathbf{Z}))=\mathbf{g}(\nabla_{\mathbf{W}}\mathbf{X},\mathbf{Z})+\mathbf{g}(\mathbf{X},\nabla_{\mathbf{W}}\mathbf{Z}),
\]
we obtain
\[
\nabla_{\mathbf{W}}\bigl(\mathbf{R}(\mathbf{X},\mathbf{Y})\mathbf{Z}\bigr)
=
-\Bigl(
\textcolor{red}{\mathbf{g}(\nabla_{\mathbf{W}}\mathbf{Y},\mathbf{Z})\mathbf{X}}
+
\textcolor{blue}{\mathbf{g}(\mathbf{Y},\nabla_{\mathbf{W}}\mathbf{Z})\mathbf{X}}
+
\textcolor{teal}{\mathbf{g}(\mathbf{Y},\mathbf{Z})\nabla_{\mathbf{W}}\mathbf{X}}
\]
\[
\quad
-
\textcolor{orange}{\mathbf{g}(\nabla_{\mathbf{W}}\mathbf{X},\mathbf{Z})\mathbf{Y}}
-
\textcolor{magenta}{\mathbf{g}(\mathbf{X},\nabla_{\mathbf{W}}\mathbf{Z})\mathbf{Y}}
-
\textcolor{olive}{\mathbf{g}(\mathbf{X},\mathbf{Z})\nabla_{\mathbf{W}}\mathbf{Y}}
\Bigr).
\]

On the other hand,
\[
\mathbf{R}(\nabla_{\mathbf{W}}\mathbf{X},\mathbf{Y})\mathbf{Z}
=
-\bigl(
\textcolor{teal}{\mathbf{g}(\mathbf{Y},\mathbf{Z})\nabla_{\mathbf{W}}\mathbf{X}}
-
\textcolor{orange}{\mathbf{g}(\nabla_{\mathbf{W}}\mathbf{X},\mathbf{Z})\mathbf{Y}}
\bigr),
\]
\[
\mathbf{R}(\mathbf{X},\nabla_{\mathbf{W}}\mathbf{Y})\mathbf{Z}
=
-\bigl(
\textcolor{red}{\mathbf{g}(\nabla_{\mathbf{W}}\mathbf{Y},\mathbf{Z})\mathbf{X}}
-
\textcolor{olive}{\mathbf{g}(\mathbf{X},\mathbf{Z})\nabla_{\mathbf{W}}\mathbf{Y}}
\bigr),
\]
\[
\mathbf{R}(\mathbf{X},\mathbf{Y})\nabla_{\mathbf{W}}\mathbf{Z}
=
-\bigl(
\textcolor{blue}{\mathbf{g}(\mathbf{Y},\nabla_{\mathbf{W}}\mathbf{Z})\mathbf{X}}
-
\textcolor{magenta}{\mathbf{g}(\mathbf{X},\nabla_{\mathbf{W}}\mathbf{Z})\mathbf{Y}}
\bigr).
\]

Substituting these expressions into the definition of $(\nabla_{\mathbf{W}}\mathbf{R})(\mathbf{X},\mathbf{Y})\mathbf{Z}$ shows that each colored term in $
\nabla_{\mathbf{W}}\bigl(\mathbf{R}(\mathbf{X},\mathbf{Y})\mathbf{Z}\bigr)
$ cancels exactly with the term of the same color arising from one of the other three summands. Therefore,
\[
(\nabla_{\mathbf{W}}\mathbf{R})(\mathbf{X},\mathbf{Y})\mathbf{Z}=0
\]
for all vector fields $\mathbf{W},\mathbf{X},\mathbf{Y},\mathbf{Z}$. Consequently,
\[
\nabla^{k}\mathbf{R}=0
\]
for every $k\geq 1$.

We also prove the global hypotheses we will use. If $c(t)=(x(t),y(t))$ is a piecewise $C^1$ curve in $\mathbb H^2$, then
\[
L_{\mathbf{g}}(c)
=\int\frac{\sqrt{|x'(t)|^2+|y'(t)|^2}}{y(t)}\,dt
\geq\int\frac{|y'(t)|}{y(t)}\,dt
\geq\bigl|\log y(b)-\log y(a)\bigr|.
\]
Let $q_j=(x_j,y_j)$ be a Cauchy sequence for the Riemannian distance. The preceding inequality shows that $(\log y_j)$ is Cauchy; in particular, for large $j$, $a\leq y_j\leq b$ for some constants $0<a<b<\infty$. Given two large indices, choose a curve between $q_j$ and $q_k$ whose length is less than $2d_{\mathbf{g}}(q_j,q_k)$. If this length is less than $1$, the inequality already proved, applied to each subarc, shows that the coordinate $y$ along the entire curve is at most $e b$. Consequently,
\[
\|x_j-x_k\|
\leq e b\,L_{\mathbf{g}}(c)
\leq 2e b\,d_{\mathbf{g}}(q_j,q_k).
\]
Thus, $(x_j)$ is also Cauchy and $q_j$ converges in the upper half-plane. Let $q=(x,y)$ be its Euclidean limit. For large $j$, the straight segment from $q_j$ to $q$ lies in a region where the second coordinate is at least $\frac{a}{2}$. In this region, for every vector $v$, $|v|_{\mathbf{g}}\leq 2a^{-1}|v|$; consequently, $d_{\mathbf{g}}(q_j,q)\leq 2a^{-1}|q_j-q|\to0$. This proves completeness. The upper half-plane is convex and hence simply connected.

The Cartan--Hadamard Theorem~\ref{teo:cartan-hadamard} now applies: for every $p\in \mathbb{H}^{2}$, the exponential map
\[
\exp_{p}\colon T_{p}\mathbb{H}^{2}\longrightarrow \mathbb{H}^{2}
\]
is a global diffeomorphism. In particular,
\[
\operatorname{inj}(p)=+\infty
\]
for every $p\in \mathbb{H}^{2}$, and therefore
\[
\operatorname{inj}(\mathbb{H}^{2})=+\infty.
\]

We conclude that $(\mathbb{H}^{2},\mathbf{g})$ has bounded geometry.

\item The unit sphere $(S^n,\mathbf{g}_{\mathrm{can}})$ has bounded geometry. Its sectional curvature is constant and equal to $1$, so Theorem~\ref{thm: tensor curvatura constante} gives
\[
\mathbf{R}(\mathbf{X},\mathbf{Y})\mathbf{Z}
=\mathbf{g}(\mathbf{Y},\mathbf{Z})\mathbf{X}
-\mathbf{g}(\mathbf{X},\mathbf{Z})\mathbf{Y}.
\]
Compatibility $\nabla \mathbf{g}=0$ implies, by the same Leibniz rule used in the hyperbolic example, that $\nabla\mathbf{R}=0$ and hence $\nabla^k\mathbf{R}=0$ for every $k\geq1$.

It remains to compute the injectivity radius, since Cartan--Hadamard does not apply here. For $p\in S^n$ and $v\in T_pS^n\setminus\{0\}$, the geodesic with initial data $(p,v)$ is
\[
\gamma_v(t)
=\cos\bigl(t|v|_{\mathbf{g}_{\mathrm{can}},p}\bigr)p
+\sin\bigl(t|v|_{\mathbf{g}_{\mathrm{can}},p}\bigr)
\frac{v}{|v|_{\mathbf{g}_{\mathrm{can}},p}};
\]
this formula is verified by differentiating twice and observing that $\gamma_v''(t)=-|v|_{\mathbf{g}_{\mathrm{can}},p}^2\gamma_v(t)$, whose acceleration is normal to the sphere. Consequently,
\[
\exp_p(v)
=\cos\bigl(|v|_{\mathbf{g}_{\mathrm{can}},p}\bigr)p
+\sin\bigl(|v|_{\mathbf{g}_{\mathrm{can}},p}\bigr)
\frac{v}{|v|_{\mathbf{g}_{\mathrm{can}},p}}.
\]
If $|v|_{\mathbf{g}_{\mathrm{can}},p}<\pi$, the Euclidean inner product with $p$ determines $|v|_{\mathbf{g}_{\mathrm{can}},p}
=\arccos\langle p,\exp_p(v)\rangle_{\mathbb R^{n+1}}$, and the orthogonal component then determines $\frac{v}{|v|_{\mathbf{g}_{\mathrm{can}},p}}$. Thus, $\exp_p$ is injective on $B_{\mathbf{g}_{\mathrm{can}},p}(0,\pi)$. We also verify that it is a local diffeomorphism on this ball. If $v=re$, with $0<r<\pi$ and $|e|=1$, write $w=be+w_\perp$, with $w_\perp\perp e$. Differentiating the explicit formula gives
\[
 d\exp_p(v)w
 =b(-\sin r\,p+\cos r\,e)+\frac{\sin r}{r}w_\perp,
 \qquad
 |d\exp_p(v)w|^2=b^2+\left(\frac{\sin r}{r}\right)^2|w_\perp|^2.
\]
The two vector summands are orthogonal, and $\sin r/r>0$. The differential is therefore invertible; at $v=0$ it is the identity. Injectivity and the inverse function theorem show that $\exp_p$ is a diffeomorphism onto its image on the indicated ball. In contrast, all vectors of $\mathbf{g}_{\mathrm{can},p}$ norm equal to $\pi$ are mapped to the antipodal point $-p$. We conclude that
\[
\operatorname{inj}(p)=\pi,
\qquad
\operatorname{inj}(S^n)=\pi.
\]
In particular, the sphere also satisfies both conditions of Definition~\ref{def: geometria acotada}.

In general, a Riemannian manifold of constant sectional curvature has a curvature tensor expressed algebraically in terms of the metric, and hence $\nabla\mathbf{R}=0$. The models $\mathbb R^n$, $S^n$, and $\mathbb H^n$ thus have bounded curvature derivatives of all orders. Their injectivity radii are $+\infty$, $\pi$, and $+\infty$, respectively, for the normalizations used here.
\end{enumerate}
\end{example}

Compact Riemannian manifolds have bounded geometry, providing many examples of such manifolds.
\begin{theorem}\label{teo:geometria-acotada-variedad-riemanniana-compacta-geometria-acotada}\index{bounded geometry@bounded geometry!compact manifolds}
Let $(M,\mathbf{g})$ be a compact Riemannian manifold without boundary. Then $(M,\mathbf{g})$ has bounded geometry.
\end{theorem}

\begin{proof}
The connected components of a manifold are open and form a cover of $M$. Compactness gives a finite subcover, so there are only finitely many components. Write
\[
M=\bigsqcup_{j=1}^{N} M_{j},
\]
where each $M_j$ is a connected component of $M$. Since connected components are closed in every topological space, each $M_j$ is closed in $M$ and hence compact.

Consequently, each Riemannian manifold $(M_j,\mathbf{g}\restriction_{M_j})$ is compact and therefore complete. For each $j\in\{1,\dots,N\}$, Proposition~\ref{prop: continuidad radio inyectividad conexa y completa} implies that the function
\[
\operatorname{inj}\colon M_j\longrightarrow (0,\infty]
\]
is continuous. (Note that $\operatorname{inj}_M(p)=\operatorname{inj}_{M_j}(p)\ \forall p \in M_j$, since $M_j$ is a connected component of $M$ and the geodesics of a manifold lie in a single connected component; this is why we simply denote it by $\operatorname{inj}(p)$.) Moreover, since $\operatorname{inj}(p)>0$ for every $p\in M_j$, there exists $p_j\in M_j$ such that
\[
\operatorname{inj}(p_j)=\min_{p\in M_j}\operatorname{inj}(p)>0.
\]
Define
\[
c:=\min_{j\in\{1,\ldots,N\}}\operatorname{inj}(p_j)>0.
\]
Then $\operatorname{inj}(p)\geq c$ for every $p\in M$, and consequently $\operatorname{inj}(M)>0$. On the other hand, for each $k\in \mathbb{N}_{0}$, the tensor $\nabla^{k}R$ is smooth on $M$, so the function
\[
p\longmapsto |\nabla^{k}R(p)|_{\mathbf{g}}
\]
is continuous. Since $M$ is compact, for each $k\in \mathbb{N}_{0}$ there exists $C_{k}>0$ such that $|\nabla^{k}R|_{\mathbf{g}}\leq C_{k}$. Therefore, $(M,\mathbf{g})$ has bounded geometry.
\end{proof}

A localized norm depends on three choices: the charts, their domains, and the cutoff functions. It is useful to collect them into a single object so that we can compare two constructions without losing track of any of these choices.

\begin{definition}[Trivialization]\label{def:geometria-acotada-trivializacion}\index{trivialization@trivialization}
Let $(M,\mathbf{g})$ be a Riemannian manifold without boundary. We say that $\mathcal{T}=(U_{i},\phi_{i},h_{i})_{i\in J}$ is a \textit{trivialization} of $M$ if $(U_{i},\phi_{i})_{i\in J}$ is an atlas of $M$ and $(h_{i})_{i\in J}$ is a partition of unity subordinate to the open cover $(U_{i})_{i\in J}$.
\end{definition}

For an infinite cover, local finiteness alone does not allow estimates to be summed with a common constant. We need a uniform bound for the number of charts that may be involved in each overlap.

\begin{definition}\label{def:geometria-acotada-espacio-topologico-cubierta-abierta-uniformemente-localmente}\index{uniformly locally finite cover}
If $X$ is a topological space and $(U_{i})_{i\in J}$ is an open cover of $X$, we say that the cover is uniformly locally finite if there exists $L>0$ such that, for each $i\in J$,
\[\#\{j\in J\mid U_{j}\cap U_{i}\neq \varnothing\}\leq L.\]
\end{definition}
We say that a trivialization is \textit{uniformly locally finite} if its underlying cover has this property.
\begin{figura}[htbp]
    \includegraphics[width=0.35\textwidth]{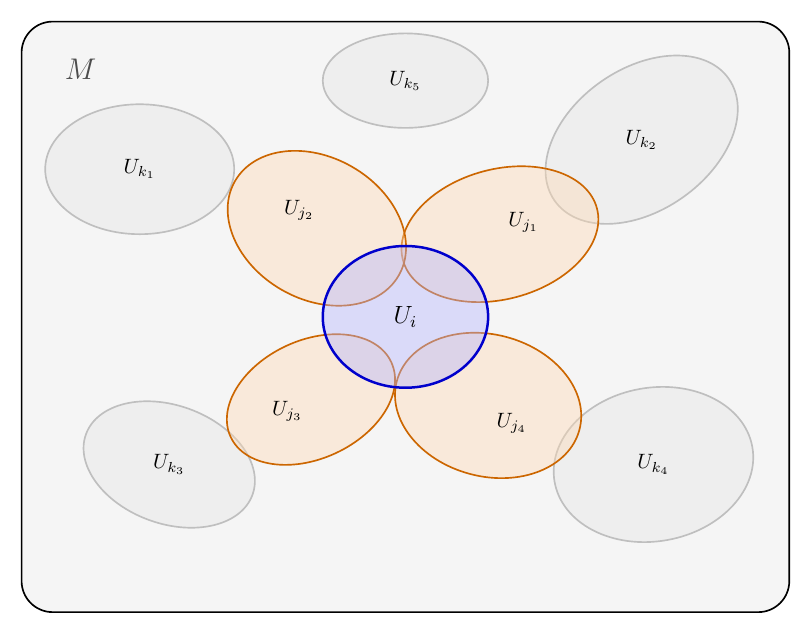}
    \caption{A uniformly locally finite cover. For any open set $U_i$ in the cover, shown in blue, the number of open sets in the cover meeting it, shown in orange, is bounded by a constant $L$ independent of $i$.}
    \label{fig:cubierta-uniformemente-localmente-finita}
\end{figura}
\begin{definition}\label{def:trivializacion-admisible}\index{admissible trivialization@admissible trivialization}
Let $(M,\mathbf{g})$ be a Riemannian manifold without boundary and with bounded geometry. If $\mathcal{T}=(U_{i},\phi_{i},h_{i})_{i\in I}$ is a uniformly locally finite trivialization of $M$, we say that $\mathcal{T}$ is \textit{admissible} if the following conditions hold:
\begin{enumerate}[label=(B\arabic*), ref=(B\arabic*)]
\item\label{item:B1}\index{admissible trivialization@admissible trivialization!uniform compatibility with geodesic atlases}
Let \(r_0>0\) be a radius for which condition (2) of Theorem~\ref{teo: equivalencias de geometria acotada} holds. There exists a uniform geodesic atlas
\[
 \mathcal A=(U_j^{\operatorname{geo}},
 \phi_j^{\operatorname{geo}})_{j\in J}
 \]
of common radius \(0<r\leq \frac{r_0}{10}\), with centers \(p_j\), for which the balls \(B_{\mathbf{g}}(p_j,\frac{r}{4})\) are pairwise disjoint and the balls \(B_{\mathbf{g}}(p_j,\frac{r}{2})\) cover \(M\). Moreover, there exists \(N\in\mathbb N\) such that
\[
 \displaystyle\sup_{i\in I}
 \#\{j\in J\mid U_i\cap U_j^{\operatorname{geo}}\neq\varnothing\}
 \leq N
 \]
and
\[
 \displaystyle\sup_{j\in J}
 \#\{i\in I\mid U_i\cap U_j^{\operatorname{geo}}\neq\varnothing\}
 \leq N.
 \]
Finally, for each \(k\in\mathbb N_0\), there exists \(C_k>0\) with the following property: if \(U_i\cap U_j^{\operatorname{geo}}\neq\varnothing\) and \(|\alpha|\leq k\), then
\[
 \displaystyle\sup_{x\in\phi_j^{\operatorname{geo}}(U_i\cap U_j^{\operatorname{geo}})}
 \left|D^{\alpha}
 \bigl(\phi_{i}\circ(\phi_{j}^{\operatorname{geo}})^{-1}\bigr)(x)\right|
 \leq C_{k},
 \]
and
\[
 \displaystyle\sup_{y\in\phi_i(U_i\cap U_j^{\operatorname{geo}})}
 \left|D^{\alpha}
 \bigl(\phi_{j}^{\operatorname{geo}}\circ\phi_i^{-1}\bigr)(y)\right|
 \leq C_k.
 \]

\item\label{item:B2}
For each $k\in \mathbb{N}_{0}$, there exists $c_{k}>0$ such that, for every $i\in I$ and every multi-index $|\alpha|\leq k$,
\[
 \displaystyle\sup_{x\in\phi_i(U_i)}
 \left|D^{\alpha}(h_{i}\circ \phi_{i}^{-1})(x)\right|
 \leq c_{k}.
 \]

\item\label{item:B3}
For each $i\in I$, there exists a function \(\chi_i\in C_c^\infty(U_i)\) such that
\[
 0\leq\chi_i\leq1,
 \qquad
 \chi_i=1
 \quad\text{on a neighborhood of }\operatorname{supp}(h_i),
 \]
and, for each $k\in\mathbb N_0$, there exists $d_k>0$, independent of $i$, for which
\[
 \displaystyle\sup_{x\in\phi_i(U_i)}
 \left|D^\alpha(\chi_i\circ\phi_i^{-1})(x)\right|
 \leq d_k,
 \qquad |\alpha|\leq k.
 \]
\end{enumerate}
\end{definition}

The functions $\chi_i$ in item \ref{item:B3} will be called \textit{auxiliary cutoff functions} of the trivialization. Their presence allows a localized expression to be extended by zero, denoted by a tilde, before a coordinate change is applied. The condition is not redundant: bounds on the derivatives of $h_i$ do not by themselves provide a uniform separation between \(\operatorname{supp}(h_i)\) and the boundary of $U_i$.

\begin{lemma}[Uniform margin for cutoff functions]
\label{lem:margen-uniforme-cortes-admisibles}
The auxiliary functions in \ref{item:B3} can be chosen so that there exists $\delta>0$, independent of $i\in I$, with
\[
 \operatorname{dist}_{\mathbb R^n}
 \bigl(\phi_i(\operatorname{supp}\chi_i),
             \mathbb R^n\setminus\phi_i(U_i)\bigr)\geq\delta.
\]
The new functions still equal one on a neighborhood of $\operatorname{supp}h_i$ and retain uniform bounds on all their derivatives.
\end{lemma}

\begin{proof}
Temporarily denote the functions in \ref{item:B3} by $\chi_i^0$. Their coordinate representations, extended by zero, are smooth on $\mathbb R^n$. The bound on their first derivatives gives a constant $D>0$ such that
\[
 \bigl\|D\widetilde{\chi_i^0\circ\phi_i^{-1}}\bigr\|_\infty
 \leq D,\qquad i\in I.
\]
Choose $\vartheta\in C^\infty(\mathbb R)$ with $0\leq\vartheta\leq1$, $\vartheta=0$ on $(-\infty,1/4]$, and $\vartheta=1$ on $[1/2,\infty)$, and set $\chi_i=\vartheta\circ\chi_i^0$. If $x\in\phi_i(\operatorname{supp}\chi_i)$, then $(\chi_i^0\circ\phi_i^{-1})(x)\geq1/4$. For $y\in\mathbb R^n\setminus\phi_i(U_i)$, the extension is zero. The mean value theorem on the Euclidean segment from $x$ to $y$ implies $1/4\leq D|x-y|$. Taking the infimum over these $y$ gives the margin $\delta=(4D)^{-1}$. The Faà di Bruno formula of Theorem~\ref{faa di bruno multivariable}, applied to the fixed function $\vartheta$, provides the higher-order bounds. Since $\chi_i^0=1$ near $\operatorname{supp}h_i$, the new function also equals one there.
\end{proof}

From now on, we choose the auxiliary functions with this margin. Thus, the constants in a zero extension or a localized coordinate change incorporate a common positive distance from the artificial boundary of the charts.

\begin{definition}[Quadratic localization system]
\label{def:sistema-cuadratico-localizacion-admisible}
\index{quadratic localization system}
Let $(M,\mathbf{g})$ be a Riemannian manifold without boundary and with bounded geometry. A \textit{quadratic localization system} is a family \[
\mathcal Q=(U_i,\phi_i,h_i)_{i\in I}
\] such that $(U_i,\phi_i)_{i\in I}$ is a uniformly locally finite atlas, $h_i\in C_c^\infty(U_i)$, \(0\leq h_i\leq1\), and
\[
\displaystyle\sum_{i\in I}h_i^2=1
\qquad\text{on }M.
\]
We say that $\mathcal Q$ is \textit{admissible} if the atlas satisfies condition \ref{item:B1}, the functions $h_i$ satisfy the estimates in \ref{item:B2}, with $h_i$ in place of the partition functions, and there is a family of auxiliary cutoff functions as in \ref{item:B3}, equal to one on a neighborhood of \(\operatorname{supp}(h_i)\).
\end{definition}

\begin{remark}
If $\mathcal Q=(U_i,\phi_i,h_i)_{i\in I}$ is an admissible quadratic localization system and $\rho_i:=h_i^2$, then $\displaystyle\sum_{i\in I}\rho_i=1$ and \[
\mathcal T_{\mathcal Q}:=(U_i,\phi_i,\rho_i)_{i\in I}
\] is an admissible trivialization. Indeed, the Leibniz rule converts the uniform derivative bounds for $h_i$ into the corresponding bounds for $h_i^2$, and the same auxiliary cutoff functions serve both families. This observation distinguishes the quadratic system used in the analysis and synthesis formulas from the genuine partition of unity.
\end{remark}

\begin{remark}[Meaning and use of quadratic localization]
\label{obs:sentido-localizacion-cuadratica}
A quadratic partition is not a second definition of the function spaces. Its usefulness lies in the fact that the same factor $h_i$ enters both the analysis and the synthesis operators. Analysis sends $u$ to the family
\[
\left(
\widetilde{(h_i u)\circ\phi_i^{-1}}
\right)_{i\in I}.
\]
For synthesis, each component is first multiplied by a cutoff function whose support is compact in the chart image, reconstructed on $U_i$, extended by zero to $M$, and multiplied by $h_i$. The precise distributional definition will be given in Proposition~\ref{prop:retraccion-interpolacion-dualidad-geometria-acotada}. The identity $\displaystyle\sum_{i\in I}h_i^2=1$ makes the composition of the two operators exactly the identity before any comparison of norms. This is Amann's convention: the cutoff system appears in \cite[Lemma~3.6]{Amann2025FunctionSpaces}, the two operators in formulas~(4.5)--(4.6), and the retraction--coretraction property in Theorem~4.2 of that reference.

The local definitions of the spaces given below use, instead, an ordinary partition \((\rho_i)\), as in \cite{traces,Schneider2021}. Lemma~\ref{lem:equivalencia-localizacion-cuadratica} proves that, when \(\rho_i=h_i^2\), both localizations satisfy the explicit inequalities with the positive constants $c_H,C_H,c_F,C_F$ of that lemma, independent of the section. It is therefore useful to distinguish the two families in the notation, even though they describe the same space.
\end{remark}

\begin{lemma}[Uniform quadratic normalization]
\label{lem:normalizacion-cuadratica-uniforme}
Let $(M,\mathbf{g})$ be a Riemannian manifold without boundary and with bounded geometry, and let \(\mathcal T=(U_i,\phi_i,\psi_i)_{i\in I}\) be an admissible trivialization such that \(0\leq\psi_i\leq1\) and \(\operatorname{supp}\psi_i\Subset U_i\) for every \(i\). Define
\[
S(x):=\left(\displaystyle\sum_{j\in I}\psi_j(x)^2\right)^{\frac{1}{2}},
\qquad
h_i(x):=\frac{\psi_i(x)}{S(x)}.
\]
Then \(\mathcal Q=(U_i,\phi_i,h_i)_{i\in I}\) is an admissible quadratic localization system.
\end{lemma}

\begin{proof}
The sum defining \(S\) is locally finite. If \(L\) is a uniform bound for the number of members of the cover that can contain a given point, the equality \(\displaystyle\sum_{j\in I}\psi_j=1\) and the Cauchy--Schwarz inequality give
\[
1
=
\left(\displaystyle\sum_{j\in I}\psi_j(x)\right)^2
\leq
L\displaystyle\sum_{j\in I}\psi_j(x)^2.
\]
Since the functions are nonnegative, we also have \(\displaystyle\sum_{j\in I}\psi_j(x)^2\leq(\displaystyle\sum_{j\in I}\psi_j(x))^2=1\). Therefore,
\[
L^{-\frac{1}{2}}\leq S(x)\leq1
\qquad\text{for every }x\in M.
\]
In particular, \(S\) is smooth, \(0\leq h_i\leq1\), \(\operatorname{supp}h_i\subseteq\operatorname{supp}\psi_i\Subset U_i\), and
\[
\displaystyle\sum_{i\in I}h_i^2
=
S^{-2}\displaystyle\sum_{i\in I}\psi_i^2
=1.
\]

It remains to verify uniformity of the derivatives. Fix a chart \((U_i,\phi_i)\). Only indices \(j\) for which \(U_i\cap U_j\neq\varnothing\) contribute, and there are at most \(L\) of them. In chart \(i\), the representation of \(\psi_j\) is understood to be extended by zero outside \(\phi_i(U_i\cap U_j)\). This extension is smooth: points where \(\psi_j\) may be nonzero lie inside \(\operatorname{supp}(\psi_j)\Subset U_j\), whereas outside that support the function vanishes on a neighborhood. On the overlap,
\[
\psi_j\circ\phi_i^{-1}
=
(\psi_j\circ\phi_j^{-1})\circ
(\phi_j\circ\phi_i^{-1}).
\]
Bounds for the coordinate transitions follow directly from \ref{item:B1}: near each point of the overlap, insert a geodesic chart from the atlas appearing in that condition and write
\[
\phi_j\circ\phi_i^{-1}
=
\bigl(\phi_j\circ(\phi_\nu^{\operatorname{geo}})^{-1}\bigr)
\circ
\bigl(\phi_\nu^{\operatorname{geo}}\circ\phi_i^{-1}\bigr).
\]
Condition \ref{item:B2}, these bounds, and the Faà di Bruno formula of Theorem~\ref{faa di bruno multivariable} bound all derivatives of these functions uniformly. The Leibniz rule and the bound \(L\) then bound all derivatives of \(\displaystyle\sum_{j\in I_i}(\psi_j\circ\phi_i^{-1})^2\).

This last function takes values in the compact interval \([L^{-1},1]\). Applying Faà di Bruno once more to \(t\mapsto t^{-\frac{1}{2}}\), whose derivatives are bounded on that interval, gives uniform bounds for the derivatives of \(S^{-1}\circ\phi_i^{-1}\). Finally, applying the Leibniz rule to \[
h_i\circ\phi_i^{-1}
=
(\psi_i\circ\phi_i^{-1})(S^{-1}\circ\phi_i^{-1})
\] gives the estimates in \ref{item:B2} for \(h_i\). This proves admissibility: the auxiliary cutoff functions of the original trivialization still equal one on a neighborhood of \(\operatorname{supp}(h_i)\subseteq\operatorname{supp}(\psi_i)\) and retain the estimates in \ref{item:B3}.
\end{proof}

\begin{remark}\label{obs:geometria-acotada-item-b1-satisface-algun-atlas-geodesico}
The coordinate-transition estimates in \ref{item:B1} do not depend on the particular normal atlas. Indeed, if \(0<r'\leq r_0\), extend the charts of radii \(r\) and \(r'\) to radius \(\displaystyle \max\{r,r'\}\leq r_0\) when necessary, apply the additional assertion of Theorem~\ref{teo: equivalencias de geometria acotada}, and restrict again.

The combinatorial condition is also independent of the atlas when the second atlas, of radius \(0<r'\leq \frac{r_0}{10}\), has disjoint balls \(B_{\mathbf{g}}(p_j',\frac{r'}{4})\) and the balls \(B_{\mathbf{g}}(p_j',\frac{r'}{2})\) cover \(M\). To verify this, if a chart of radius \(r\) meets a chart of radius \(r'\), the disjoint ball of radius \(\frac{r'}{4}\) associated with the second lies in a ball whose radius does not exceed \(r+r'+\frac{r'}{4}<r_0\). The uniform bounds on the metric and its inverse in normal charts, given by Theorem~\ref{teo: equivalencias de geometria acotada}, place the eigenvalues of its matrix between two positive constants independent of the center. Integrating the square root of its determinant gives, for each fixed radius less than $r_0$, uniform upper and lower bounds for the volume of the balls. Comparing the large ball with the disjoint balls of radius $\frac{r'}{4}$ bounds their number uniformly. Interchanging \(r\) and \(r'\) gives the other overlap bound. Combining this fact with the two overlap bounds in \ref{item:B1} also transfers the cross-overlap bounds between \(\mathcal T\) and the second atlas.

Moreover, \ref{item:B1} controls transitions between two admissible trivializations. Let \(\mathcal T=(U_i,\phi_i)_{i\in I}\) and \(\mathcal T^\sharp
   =(U_j^\sharp,\phi_j^\sharp)_{j\in J}\). By the preceding two paragraphs, we may fix a common uniform geodesic atlas
\((U_\nu^{\operatorname{geo}},
   \phi_\nu^{\operatorname{geo}})_{\nu\in K}\)
relative to which both satisfy the bounds in \ref{item:B1}. If \(p\in U_i\cap U_j^\sharp\), choose \(\nu\) with \(p\in U_\nu^{\operatorname{geo}}\). On a neighborhood of \(p\) contained in the triple overlap, we have exactly
\[
 \phi_j^\sharp\circ\phi_i^{-1}
 =
 \bigl(\phi_j^\sharp\circ
       (\phi_\nu^{\operatorname{geo}})^{-1}\bigr)
 \circ
 \bigl(\phi_\nu^{\operatorname{geo}}\circ\phi_i^{-1}\bigr).
 \]
The Faà di Bruno formula of Theorem~\ref{faa di bruno multivariable} expresses every derivative of order \(k\) of the left-hand side at \(p\) as a finite sum of products of derivatives of the two factors of orders at most \(k\). The bounds in \ref{item:B1} are independent of \(i,j,\nu\); since \(p\) is arbitrary, for each \(k\in\mathbb N_0\) there exists \(C_k\), independent of \(i,j\), such that
\[
 \left\|D^\alpha\bigl(\phi_j^\sharp\circ\phi_i^{-1}\bigr)\right\|_{L^\infty(\phi_i(U_i\cap U_j^\sharp))}
 \leq C_k,
 \qquad |\alpha|\leq k.
 \]
Interchanging the two trivializations gives the same estimate for the inverse transition. The case \(\mathcal T^\sharp=\mathcal T\) provides, in particular, uniform control of the internal transitions of an admissible trivialization.
\end{remark}
Normal geodesic coordinates can be used to define a trivialization of any manifold of bounded geometry that is also \textit{uniformly locally finite}; we will call it a \textit{\textbf{geodesic trivialization}}.
\begin{lemma}[Uniform control of geodesic ball volumes]
\label{lema: volumen bolas geometria acotada}\index{uniform control of geodesic ball volumes@uniform control of geodesic ball volumes}
Let $(M,\mathbf{g})$ be a Riemannian manifold without boundary and with bounded geometry. If \(r_0>0\) is a radius for which condition (2) of Theorem~\ref{teo: equivalencias de geometria acotada} holds, then, for each \(0<r\leq r_0\), there exist constants \(v_{-}(r),v_{+}(r)>0\) such that
\[
v_{-}(r)
\le
\lambda_{\mathbf{g}}(B_{\mathbf{g}}(p,r))
\le
v_{+}(r),
\qquad \forall p\in M.
\]
\end{lemma}

\begin{proof}
Fix \(0<r\leq r_0\) and let \(p\in M\). The Riemannian exponential map at \(p\) induces a normal chart on \(B_{\mathbf{g}}(p,r)\). Fix a linear isomorphism \[
B_p\colon \mathbb{R}^{n}\longrightarrow T_pM
\] taking the standard basis of $\mathbb{R}^{n}$ to an orthonormal basis of $T_pM$, and consider the chart
\[
\phi_p:=B_p^{-1}\circ \exp_p^{-1}\colon B_{\mathbf{g}}(p,r)\longrightarrow B_{\mathrm{euc}}(0,r).
\]

Denote by $g_{ij}^{(p)}$ and $g_{(p)}^{ij}$ the components of the metric and its inverse in these coordinates. Since $(M,\mathbf{g})$ has bounded geometry, Theorem~\ref{teo: equivalencias de geometria acotada} implies that there exists a constant $C>0$, independent of $p$, such that
\[
|g_{ij}^{(p)}(q)|\le C
\qquad\text{and}\qquad
|g_{(p)}^{ij}(q)|\le C,
\qquad \forall q\in B_{\mathbf{g}}(p,r),
\]
for all $i,j\in \{1,\dots,n\}$.

Fix $q\in B_{\mathbf g}(p,r)$ and let $\lambda_1(q),\dots,\lambda_n(q)$ be the eigenvalues of the coefficient matrix $(g_{ij}^{(p)}(q))_{i,j=1}^n$. This matrix is symmetric and positive definite. The bound on its entries implies that its Frobenius norm is at most $nC$; therefore,
\[
 \lambda_\ell(q)\leq nC,
 \qquad \ell=1,\dots,n.
\]
The inverse matrix has coefficients $(g_{(p)}^{ij}(q))$ and its Frobenius norm is also at most $nC$. Its eigenvalues are $\lambda_\ell(q)^{-1}$, whence
\[
 \lambda_\ell(q)^{-1}\leq nC
 \quad\Longrightarrow\quad
 \lambda_\ell(q)\geq(nC)^{-1}.
\]
Multiplying the $n$ inequalities gives directly
\[
 (nC)^{-n}
 \leq \det\bigl(\mathbf{g}^{(p)}(q)\bigr)
 \leq (nC)^n,
 \qquad q\in B_{\mathbf g}(p,r),
\]
uniformly in $p$. Consequently,
\[
 (nC)^{-\frac n2}
 \leq \sqrt{\det\bigl(\mathbf{g}^{(p)}(q)\bigr)}
 \leq (nC)^{\frac n2}.
\]

Using the local expression for the Riemannian measure, we obtain
\[
\lambda_{\mathbf{g}}(B_{\mathbf{g}}(p,r))
=
\int_{B_{\mathrm{euc}}(0,r)}
\sqrt{\det\bigl(\mathbf{g}^{(p)}(\phi_{p}^{-1}(x))\bigr)}\,dx.
\]
By the preceding bounds,
\[
(nC)^{-\frac{n}{2}}\,\lambda_{n}(B_{\mathrm{euc}}(0,r))
\leq
\lambda_{\mathbf{g}}(B_{\mathbf{g}}(p,r))
\le
(nC)^{\frac{n}{2}}\,\lambda_{n}(B_{\mathrm{euc}}(0,r)).
\]

It therefore suffices to define
\[
v_-(r):=(nC)^{-\frac{n}{2}}\lambda_{n}(B_{\mathrm{euc}}(0,r))
\qquad\text{and}\qquad
v_+(r):=(nC)^{\frac{n}{2}}\lambda_{n}(B_{\mathrm{euc}}(0,r)).
\]
This concludes the proof.
\end{proof}
\begin{theorem}[Geodesic trivialization]
\label{thm:trivializacion-geodesica}\index{geodesic trivialization@geodesic trivialization}
Let $(M,\mathbf{g})$ be a Riemannian manifold without boundary and with bounded geometry. Then there exist $r>0$ and an admissible trivialization \[
\mathcal{T}^{\operatorname{geo}}
=
\bigl(B_{\mathbf{g}}(p_i,2r),\phi_{i},h_{i}^{\operatorname{geo}}\bigr)_{i\in J}
\] of $M$ such that:
\begin{enumerate}
\item For each $i\in J$,
\[
 \phi_{i}=B_{p_{i}}^{-1}\circ \exp_{p_{i}}^{-1}\colon B_{\mathbf{g}}(p_i,2r)\longrightarrow B_{\mathrm{euc}}(0,2r)\subseteq \mathbb{R}^{n}
 \]
is a normal coordinate chart, where $B_{p_{i}}\colon \mathbb{R}^n\longrightarrow T_{p_{i}}M$ is an isomorphism.

\item The balls $B_{\mathbf{g}}(p_i,\displaystyle\frac{r}{2})$ are pairwise disjoint.

\item $M=\displaystyle\bigcup_{i\in J} B_{\mathbf{g}}(p_i,r)$, and the cover $(B_{\mathbf{g}}(p_i,2r))_{i\in J}$ is uniformly locally finite.

\item $(h_{i}^{\operatorname{geo}})_{i\in J}$ is a partition of unity subordinate to $(B_{\mathbf{g}}(p_i,2r))_{i\in J}$.

\item For every $k\in\mathbb{N}_{0}$, there exists $C_{k}>0$ such that
\[
 \left|D^\alpha(h_{j}^{\operatorname{geo}}\circ \phi_{i}^{-1})\right|\le C_k,
 \qquad |\alpha|\le k,
 \]
whenever $B_{\mathbf{g}}(p_i,2r)\cap B_{\mathbf{g}}(p_j,2r)\neq\varnothing$.

\item There exist functions \(\chi_i^{\operatorname{geo}}\in C_c^\infty(B_{\mathbf{g}}(p_i,2r))\), equal to one on a neighborhood of \(\operatorname{supp}(h_i^{\operatorname{geo}})\), whose representations in the normal charts have derivatives of every order uniformly bounded in $i$.
\end{enumerate}
\end{theorem}

\begin{proof}
Since $(M,\mathbf{g})$ has bounded geometry, \(\operatorname{inj}(M)>0\). Let \(r_0>0\) be a radius for which condition (2) of Theorem~\ref{teo: equivalencias de geometria acotada} holds. Choose
\[
0<r<
\min\left\{
\frac{\operatorname{inj}(M)}{5},
\frac{r_0}{20}
\right\}.
\]

Consider
\[
\mathcal{F}
=
\{A\subseteq M \mid d(p,q)\ge r \text{ for all }p,q\in A,\text{ with } p\neq q\}.
\]
We can order $\mathcal{F}$ by set inclusion $\subseteq$, obtaining a partial order. If $\mathcal{C}\subseteq \mathcal{F}$ is a chain, then $A_{\mathcal{C}}=\displaystyle\bigcup_{A\in\mathcal{C}}A$ belongs to $\mathcal{F}$, since, given $p,q\in A_{\mathcal{C}}$, total ordering of $(\mathcal{C},\subseteq )$ ensures that both belong to the same member of the chain. Thus, $A_{\mathcal{C}}$ is an upper bound for $\mathcal{C}$, so Zorn's Lemma, Theorem~\ref{teo:lema-zorn}, guarantees the existence of a maximal set $A=\{p_i\}_{i\in J}$. Maximality preserves the separation $d(p_i,p_j)\geq r$ whenever $i\neq j$, which implies that the balls $B_{\mathbf{g}}(p_i,\displaystyle\frac{r}{2})$ are disjoint.

We claim that $M=\displaystyle\bigcup_{i\in J} B_{\mathbf{g}}(p_i,r)$. Otherwise, there would exist $x\in M$ with $d(x,p_i)\ge r$ for every $i\in J$. In particular, this gives $x\notin A$, so
\[
A\subsetneq A\cup\{x\}\in \mathcal{F},
\]
contradicting maximality of $A$.

For each $i\in J$, since $2r<\operatorname{inj}(M)$, the Riemannian exponential map at $p_i$ induces a normal chart on $B_{\mathbf{g}}(p_i,2r)$. Define
\[
\phi_i:=B_{p_{i}}^{-1}\circ \exp_{p_i}^{-1}\colon B_{\mathbf{g}}(p_i,2r)\longrightarrow B_{\mathrm{euc}}(0,2r),
\]
where $B_{p_{i}}\colon \mathbb{R}^{n}\longrightarrow T_{p_{i}}M$ is an isomorphism taking the standard basis to an orthonormal basis. This proves (a).

We now prove that the cover $(B_{\mathbf{g}}(p_i,2r))_{i\in J}$ is uniformly locally finite. For fixed $i\in J$, define
\[
J_i:=\{j\in J \mid B_{\mathbf{g}}(p_i,2r)\cap B_{\mathbf{g}}(p_j,2r)\neq\varnothing\}.
\]
If $j\in J_i$, then $B_{\mathbf{g}}(p_i,2r)\cap B_{\mathbf{g}}(p_j,2r)\neq\varnothing$, whence
\[
d(p_i,p_j)<4r,
\qquad\text{and therefore}\qquad
p_j\in B_{\mathbf{g}}(p_i,4r).
\]
The balls $B_{\mathbf{g}}(p_j,\displaystyle\frac{r}{2})$ are pairwise disjoint, and moreover,
\[
B_{\mathbf{g}}(p_j,\displaystyle\frac{r}{2})\subseteq B_{\mathbf{g}}(p_i,\displaystyle\frac{9r}{2}),
\qquad \forall j\in J_i.
\]
By Lemma~\ref{lema: volumen bolas geometria acotada}, there exist constants $v_-,v_+>0$ such that
\[
\lambda_{\mathbf{g}}\bigl(B_{\mathbf{g}}(p,\displaystyle\frac{r}{2})\bigr)\ge v_-,
\qquad
\lambda_{\mathbf{g}}\bigl(B_{\mathbf{g}}(p,\displaystyle\frac{9r}{2})\bigr)\le v_+,
\qquad \forall p\in M.
\]
Summing volumes over $j\in J_i$ gives
\[
\displaystyle\sum_{j\in J_i} \lambda_{\mathbf{g}}\bigl(B_{\mathbf{g}}(p_j,\displaystyle\frac{r}{2})\bigr)
\le
\lambda_{\mathbf{g}}\bigl(B_{\mathbf{g}}(p_i,\displaystyle\frac{9r}{2})\bigr),
\]
and therefore
\[
\#J_i\,v_- \le v_+,
\qquad\text{whence}\qquad
\#J_i\le \frac{v_+}{v_-}.
\]
This implies that the cover $(B_{\mathbf{g}}(p_i,2r))_{i\in J}$ is uniformly locally finite, proving (c).

We now construct a partition of unity. Use the flat function from the proof of Theorem~\ref{particionesdelaunidad}, which we write again here:
\[
\vartheta(t):=
\begin{cases}
e^{-1/t},&t>0,\\
0,&t\leq0.
\end{cases}
\]
All its derivatives at the origin vanish. Define explicitly
\[
\psi(x):=
\frac{\vartheta\!\left((\frac{4r}{3})^2-\|x\|^2\right)}
{\vartheta\!\left((\frac{4r}{3})^2-\|x\|^2\right)
 +\vartheta\!\left(\|x\|^2-r^2\right)}.
\]
The denominator does not vanish, and flatness of $\vartheta$ at the origin shows that $\psi\in C_c^\infty(\mathbb R^n)$. Moreover,
\[
0\le\psi\le1,\qquad
\psi\equiv1 \text{ on } B_{\mathrm{euc}}(0,r),
\qquad
\operatorname{supp}(\psi)\subseteq B_{\mathrm{euc}}(0,\displaystyle\frac{3r}{2}).
\]
For each $i\in J$, define
\[
\psi_i(x):=
\begin{cases}
\psi(\phi_i(x)),& x\in B_{\mathbf{g}}(p_i,2r),\\
0,& \text{otherwise}.
\end{cases}
\]
Then $\psi_i$ is smooth and
\[
\operatorname{supp}(\psi_i)\subseteq B_{\mathbf{g}}(p_i,\displaystyle\frac{3r}{2})\subset B_{\mathbf{g}}(p_i,2r),
\]
so its extension by zero is smooth.

Define
\[
\Psi(x):=\displaystyle\sum_{i\in J} \psi_i(x).
\]
Since the cover $(B_{\mathbf{g}}(p_i,2r))_{i\in J}$ is uniformly locally finite, for every $x\in M$ there exists a neighborhood $U$ of $x$ such that the set
\[
\{i\in J \mid U\cap B_{\mathbf{g}}(p_i,2r)\neq\varnothing\}
\]
is finite. In particular, for every $y\in U$, the sum defining $\Psi(y)$ has finitely many nonzero summands, with a bound independent of $y\in U$. It follows that $\Psi$ is a smooth function on $M$.

Moreover, since \[
M=\bigcup_{i\in J} B_{\mathbf{g}}(p_i,r)
\] and $\psi_i\equiv1$ on $B_{\mathbf{g}}(p_i,r)$, for every $x\in M$ there exists $i\in J$ such that $\psi_i(x)=1$. Consequently,
\[
\Psi(x)\ge 1,
\qquad \forall x\in M.
\]
In particular, $\Psi$ does not vanish on $M$, so the functions
\[
h_i^{\operatorname{geo}}(x):=\frac{\psi_i(x)}{\Psi(x)}
\]
are well defined and smooth. Moreover, since $0\le \psi_i\le 1$ and $\Psi\ge 1$,
\[
0\le h_i^{\operatorname{geo}}\le 1,
\qquad \forall i\in J.
\]
Thus, $(h_i^{\operatorname{geo}})_{i\in J}$ is a partition of unity subordinate to $(B_{\mathbf{g}}(p_i,2r))_{i\in J}$, proving (d).

We prove the uniform derivative bounds. For fixed $i,j\in J$ such that $B_{\mathbf{g}}(p_i,2r)\cap B_{\mathbf{g}}(p_j,2r)\neq\varnothing$, in normal coordinates we have
\[
h_j^{\operatorname{geo}}\circ \phi_i^{-1}
=
\frac{\psi_j\circ \phi_i^{-1}}{\displaystyle\sum_{k\in J_i} \psi_k\circ \phi_i^{-1}},
\]
where $J_i$ is finite and has uniformly bounded cardinality.

First, we show that the functions $\psi_k\circ \phi_i^{-1}$ have uniformly bounded derivatives of every order. Fix $m\in\mathbb{N}_0$. Since $\psi$ is a fixed cutoff function on $\mathbb{R}^n$ with compact support contained in $B_{\mathrm{euc}}(0,\displaystyle\frac{3r}{2})$, there exists a constant $A_m>0$ such that
\[
|D^\alpha\psi(x)|\le A_m,
\qquad \forall x\in \mathbb{R}^n,\quad |\alpha|\le m.
\]
On the other hand, for each $k\in J_i$,
\[
\psi_k\circ \phi_i^{-1}
=
\psi\circ \phi_k\circ \phi_i^{-1}.
\]
The equality is used on \(\phi_i(B_{\mathbf{g}}(p_i,2r)\cap B_{\mathbf{g}}(p_k,2r))\). Outside this overlap, \(\psi_k\circ\phi_i^{-1}=0\). Indeed, \(\operatorname{supp}\psi\subset B_{\mathrm{euc}}(0,\frac{3r}{2})\), so \(\psi_k\) vanishes on a neighborhood of the boundary of \(B_{\mathbf{g}}(p_k,2r)\). Thus, the coordinate expression is smooth throughout $\phi_i(B_{\mathbf{g}}(p_i,2r))$, and its derivatives vanish outside the overlap. Since the geodesic coordinate transitions $\phi_k\circ \phi_i^{-1}$ have uniformly bounded derivatives through order $m$ by Theorem~\ref{teo: equivalencias de geometria acotada}, the Faà di Bruno formula of Theorem~\ref{faa di bruno multivariable} implies that there exists a constant $C_m'>0$, independent of $i$ and $k$, such that
\[
|D^\alpha(\psi_k\circ \phi_i^{-1})|\le C_m',
\qquad |\alpha|\le m.
\]

Moreover, since the sum in the denominator involves a uniformly bounded number of terms, for each $m\in\mathbb{N}_0$ there exists a constant $C_m''>0$ such that
\[
\left|D^\alpha\left(\displaystyle\sum_{k\in J_i}\psi_k\circ \phi_i^{-1}\right)\right|
\le C_m'',
\qquad |\alpha|\le m.
\]
On the other hand, for every $x\in B_{\mathrm{euc}}(0,2r)$,
\[
\displaystyle\sum_{k\in J_i}\psi_k\circ \phi_i^{-1}(x)
=
\Psi\bigl(\phi_i^{-1}(x)\bigr)\ge 1.
\]
Thus, the denominator is uniformly bounded below by $1$. Set
\[
F_i:=\displaystyle\sum_{k\in J_i}\psi_k\circ\phi_i^{-1}.
\]
The lower bound gives $\|F_i^{-1}\|_{L^\infty(B_{\mathrm{euc}}(0,2r))}\leq1$. For $|\alpha|\geq1$, applying the Leibniz rule to $F_i^{-1}F_i=1$ and isolating the term $D^\alpha(F_i^{-1})F_i$ gives the recurrence
\[
D^\alpha(F_i^{-1})
=-F_i^{-1}
\displaystyle\sum_{\substack{\beta\in\mathbb N_0^n\\\beta\leq\alpha,\ \beta\neq\alpha}}
\binom{\alpha}{\beta}
D^\beta(F_i^{-1})D^{\alpha-\beta}F_i.
\]
For $|\alpha|=1$, the right-hand side contains only $F_i^{-1}$ and first derivatives of $F_i$, which are already bounded. If the derivatives of $F_i^{-1}$ of order less than $|\alpha|$ are uniformly bounded, every term in the sum is bounded by the induction hypothesis and the bounds on $F_i$. Thus, order by order, we obtain uniform bounds for all derivatives of $F_i^{-1}$. Writing
\[
h_j^{\operatorname{geo}}\circ \phi_i^{-1}
=
(\psi_j\circ \phi_i^{-1})\cdot
F_i^{-1},
\]
the multi-index Leibniz rule and the bounds just proved imply that, for each $m\in\mathbb{N}_0$, there exists a constant $C_m>0$ such that
\[
|D^\alpha(h_j^{\operatorname{geo}}\circ \phi_i^{-1})|\le C_m,
\qquad |\alpha|\le m.
\]
This proves (e). The underlying atlas has common radius \(2r\); the balls of radius \(\frac{2r}{4}=\frac{r}{2}\) are disjoint, and those of radius \(\frac{2r}{2}=r\) cover \(M\). Moreover, \(2r<\frac{r_0}{10}\) by the choice of radius. The overlap bounds already proved and the additional assertion of Theorem~\ref{teo: equivalencias de geometria acotada} verify \ref{item:B1}, while (e) verifies \ref{item:B2}.

To construct the auxiliary cutoff functions, use the same flat function and set
\[
\chi(x):=
\frac{\vartheta\!\left((\frac{15r}{8})^2-\|x\|^2\right)}
{\vartheta\!\left((\frac{15r}{8})^2-\|x\|^2\right)
 +\vartheta\!\left(\|x\|^2-(\frac{7r}{4})^2\right)}.
\]
Then $\chi\in C_c^\infty(B_{\mathrm{euc}}(0,2r))$, $0\leq\chi\leq1$, and $\chi=1$ on $B_{\mathrm{euc}}(0,\frac{7r}{4})$. Define
\[
\chi_i^{\operatorname{geo}}
:=
\begin{cases}
\chi\circ\phi_i,&\text{on }B_{\mathbf{g}}(p_i,2r),\\
0,&\text{outside }B_{\mathbf{g}}(p_i,2r).
\end{cases}
\]
The construction of $h_i^{\operatorname{geo}}$ shows that its support is contained in $B_{\mathbf{g}}(p_i,\frac{3r}{2})$; therefore, \(\chi_i^{\operatorname{geo}}=1\) on a neighborhood of that support. In its own chart, its representation is the fixed function \(\chi\), so the required bounds are independent of $i$. This verifies \ref{item:B3} and the last item of the statement. Therefore, \(\mathcal T^{\operatorname{geo}}\) is admissible.
\end{proof}
We next present a useful lemma:
\begin{lemma}\label{lema: equivalencia uniforme metricas}
Let $(M,\mathbf{g})$ be a Riemannian manifold without boundary, and let $\mathcal{T}=(U_i,\phi_i,h_i)_{i\in I}$ be a trivialization of $M$. Suppose that there exists a constant $C>0$ such that, for every $i\in I$, the components of the metric and its inverse in the chart $(U_i,\phi_i)$ satisfy
\[
|g_{ab}^{(i)}(p)|\le C
\qquad\text{and}\qquad
|g_{(i)}^{ab}(p)|\le C
\]
for every $p\in U_i$ and all $a,b\in\{1,\dots,n\}$.

Then there exist positive constants $c_0,C_0$, independent of $i$ and $p$, such that
\[
c_0\|d(\phi_i)_p(v)\|^2
\le
|v|_{\mathbf{g}}^2
\le
C_0\|d(\phi_i)_p(v)\|^2,
\qquad \forall p\in U_i,\ \forall v\in T_pM.
\]
\end{lemma}

\begin{proof}
Fix $i\in I$ and $p\in U_i$. Let \(v\in T_pM\) and define
\[
w:=d(\phi_i)_p(v)\in \mathbb{R}^n.
\]
Since \((U_i,\phi_i)\) is a chart, the vector $v$ can be written uniquely as
\[
v=w^a\boldsymbol{\partial}_{a}\bigg|_p,
\]
so, by definition of the metric coefficients,
\begin{equation}
\label{eq:metrica-coordenada-sin-matriz-auxiliar}
 |v|_{\mathbf g}^2
 =\mathbf g_p(v,v)
 =\displaystyle\sum_{a,b=1}^n g_{ab}^{(i)}(p)w^aw^b.
\end{equation}

Let $\lambda_1(p),\dots,\lambda_n(p)$ be the eigenvalues of the symmetric positive definite matrix $(g_{ab}^{(i)}(p))_{a,b=1}^n$. The bound on its entries gives
\[
 \left(\displaystyle\sum_{a,b=1}^n|g_{ab}^{(i)}(p)|^2\right)^{1/2}
 \leq nC,
\]
and Theorem~\ref{valor propio y frobenius} implies $\lambda_\ell(p)\leq nC$. The inverse matrix has entries $g_{(i)}^{ab}(p)$, so the same argument applied to it gives
\[
 \lambda_\ell(p)^{-1}\leq nC
 \quad\Longrightarrow\quad
 \lambda_\ell(p)\geq\frac1{nC}.
\]
The spectral theorem provides an orthonormal basis $(\varepsilon_\ell)_{\ell=1}^n$ in which, if $\displaystyle w=\displaystyle\sum_{\ell=1}^n w_\ell\varepsilon_\ell$, the last member of \eqref{eq:metrica-coordenada-sin-matriz-auxiliar} is $\displaystyle \sum_{\ell=1}^n\lambda_\ell(p)w_\ell^2$. Therefore,
\[
\frac{1}{nC}\|d(\phi_i)_p(v)\|^2
\le
|v|_{\mathbf{g}}^2
\le
nC\,\|d(\phi_i)_p(v)\|^2.
\]
The constants $c_0=(nC)^{-1}$ and $C_0=nC$ are independent of $i$ and $p$, as required.
\end{proof}

\begin{remark}[Meaning of uniform comparison in coordinates]
\label{obs:comparacion-metrica-coordenada-uniforme}
If $\overline{\mathbf g}$ is the Euclidean metric and $\boldsymbol\delta_i:=\phi_i^*\overline{\mathbf g}$, then
\[
 (\boldsymbol\delta_i)_p(v,v)=\|d(\phi_i)_p(v)\|^2.
\]
Thus, the conclusion of the preceding lemma can be abbreviated as
$c_0\boldsymbol\delta_i\leq\mathbf g\leq C_0\boldsymbol\delta_i$
on the domain of each chart, in the sense of Definition~\ref{def:orden-formas-metricas}. The lemma's comparison between intrinsic and coordinate quantities consists precisely of evaluating these two metrics on the same vector $v\in T_pM$. Proposition~\ref{prop:comparacion-metricas-duales-tensores-volumen} then provides uniform bounds for covectors, tensors of any fixed type, and volume elements.

Here $c_0$ and $C_0$ are independent of the index $i$ because the hypothesis controls the coefficients in all charts and their inverses by a single constant. Applying the comparison separately on each compact set would generally give constants $c_i,C_i$; an infinite family of charts does not allow them to be replaced by common constants without also proving $\displaystyle\inf_{i\in I}c_i>0$ and $\displaystyle\sup_{i\in I} C_i<\infty$.
\end{remark}

\section{Uniform regularity as a coordinate formulation of bounded geometry}
\label{sec:regularidad-uniforme-amann-geometria-acotada}

Bounded geometry is the intrinsic notion adopted in this book. Amann's uniform regularity expresses the same control through an atlas with normalized domains; its usefulness lies in making visible the constants associated with the coefficients, cutoff functions, and localization operators. It does not introduce a uniformity on a set of configurations; the notion based on neighborhoods of the diagonal is developed in Section~\ref{sec:espacios-uniformes-eichhorn}. The spaces $BC^j$ and $BUC^j$, including the convention on a half-box, were defined in Definition~\ref{def:espacios-bc-buc}.

Set
\[
 Q^m:=(-1,1)^m,
 \qquad
 Q_+^m:=[0,1)\times(-1,1)^{m-1}.
\]
If $0<\rho<1$, write $\rho Q^m=(-\rho,\rho)^m$ and $\rho Q_+^m=[0,\rho)\times(-\rho,\rho)^{m-1}$.

\begin{definition}[Uniformly regular atlas]
\label{def:atlas-uniformemente-regular-amann}\glsadd{uniformemente-regular}
\index{normalized atlas}
\index{shrinkable atlas}
\index{uniformly regular atlas}
\index{multiplicity of a cover}
Let $M$ be a separable smooth manifold of dimension $m$, with or without boundary. An atlas $\mathfrak K=(U_\kappa,\kappa)_{\kappa\in K}$ is \textit{normalized} if \[
 \kappa\colon U_\kappa\longrightarrow Q_\kappa,
 \qquad
 Q_\kappa\in\{Q^m,Q_+^m\},
\], with $Q_+^m$ used precisely for charts meeting the boundary. It is \textit{shrinkable} if there exists a number $\rho\in(0,1)$, common to all charts, such that $\{\kappa^{-1}(\rho Q_\kappa)\mid\kappa\in K\}$ covers $M$. Its multiplicity is
\[
 \operatorname{mult}(\mathfrak K)
 :=\displaystyle\sup_{p\in M}\#\{\kappa\in K\mid p\in U_\kappa\}.
\]
The atlas is \textit{uniformly regular} if it is normalized, shrinkable, has finite multiplicity, and, for each $j\in\mathbb N_0$, there exists $C_{\mathrm{tr}}(j)>0$, independent of the charts, such that
\[
 \|\widetilde\kappa\circ\kappa^{-1}\|_{
 BUC^j(\kappa(U_\kappa\cap U_{\widetilde\kappa}),\mathbb R^m)}
 \leq C_{\mathrm{tr}}(j)
\]
if $U_\kappa\cap U_{\widetilde\kappa}\neq\varnothing$. Equivalently, every transition belongs to $BUC^\infty(\kappa(U_\kappa\cap U_{\widetilde\kappa}),\mathbb R^m)$ and hence to $BC^\infty$ on that domain, with bounds uniform over the charts for each finite-order seminorm.
\end{definition}

\begin{definition}[Uniformly regular structure]
\label{def:estructura-uniformemente-regular-amann}
\index{uniformly regular structure}
Let $M$ be a smooth manifold with or without boundary. Two uniformly regular atlases of $M$ are equivalent if there exists $N\in\mathbb N$ bounding, in both directions, the number of charts of one atlas meeting any chart of the other, and if, for each $j\in\mathbb N_0$, the mixed transitions and their inverses have $BUC^j$ norm bounded by a constant independent of the charts. An equivalence class is a \textit{uniformly regular structure}.
\end{definition}

\begin{definition}[Uniformly regular Riemannian manifold]
\label{def:riemanniana-uniformemente-regular-amann}
\index{Riemannian manifold!uniformly regular}
Let $(M,\mathbf{g})$ be a Riemannian manifold with or without boundary, and let $\mathfrak K$ be a uniformly regular atlas. For each chart, define, without ambiguity concerning pushforwards,
\[
 \mathbf{g}_\kappa:=(\kappa^{-1})^*\mathbf{g}.
\]
The metric is \textit{uniformly regular relative to $\mathfrak K$} if there exist $c_{\mathbf{g}},C_{\mathbf{g}}>0$, independent of $\kappa$, such that
\begin{equation}
\label{eq:amann-comparacion-metrica-explicita}
 c_{\mathbf{g}}|\xi|^2
 \leq \mathbf{g}_\kappa(x)(\xi,\xi)
 \leq C_{\mathbf{g}}|\xi|^2
\end{equation}
for every $x\in Q_\kappa$ and $\xi\in\mathbb R^m$, and if, for each $j\in\mathbb N_0$, there exists $G_j>0$, independent of $\kappa$, such that
\begin{equation}
\label{eq:amann-buc-metrica}
 \|(g_\kappa)_{ab}\|_{BUC^j(Q_\kappa)}\leq G_j,
 \qquad 1\leq a,b\leq m.
\end{equation}
A manifold equipped with the corresponding class is called a \textit{uniformly regular Riemannian manifold}.
\end{definition}

Bounds on the inverse metric are not an additional hypothesis. Indeed, \eqref{eq:amann-comparacion-metrica-explicita} implies $\|\mathbf{g}_\kappa^{-1}(x)\|_{\mathrm{op}}
\leq c_{\mathbf{g}}^{-1}$. Differentiating $\mathbf{g}_\kappa^{-1}\mathbf{g}_\kappa=I$ gives
\[
 D(\mathbf{g}_\kappa^{-1})[v]
 =-\mathbf{g}_\kappa^{-1}D\mathbf{g}_\kappa[v]\mathbf{g}_\kappa^{-1}.
\]
For a multi-index $\alpha\neq0$, applying the Leibniz rule to $\mathbf g_\kappa^{-1}\mathbf g_\kappa=I$ gives
\[
 0=\displaystyle\sum_{\substack{\beta\in\mathbb N_0^m\\\beta\leq\alpha}}\binom{\alpha}{\beta}
 (D^\beta\mathbf g_\kappa^{-1})
 (D^{\alpha-\beta}\mathbf g_\kappa).
\]
Separating the term $\beta=\alpha$ and multiplying on the right by $\mathbf g_\kappa^{-1}$ gives the exact recurrence
\begin{equation}
\label{eq:recurrencia-derivadas-metrica-inversa-amann}
D^\alpha\mathbf g_\kappa^{-1}
=-
\displaystyle\sum_{\substack{\beta\in\mathbb N_0^m\\\beta\leq\alpha,\ \beta\neq\alpha}}
\binom{\alpha}{\beta}
(D^\beta\mathbf g_\kappa^{-1})
(D^{\alpha-\beta}\mathbf g_\kappa)
\mathbf g_\kappa^{-1}.
\end{equation}
Order zero is bounded by $c_{\mathbf g}^{-1}$. Suppose all derivatives of the inverse of order less than $j$ are bounded, and let $|\alpha|=j$. In each summand of \eqref{eq:recurrencia-derivadas-metrica-inversa-amann}, we have $|\beta|<j$ and $1\leq|\alpha-\beta|\leq j$; the induction hypothesis, the bound $G_j$ for the metric, and the order-zero bound for its inverse control the product. The sum is finite, and its coefficients depend only on $\alpha$. This proves the bound of order $j$ and completes the induction. Thus, for each $j$, there exists a constant
\[
 G'_j=P_j(c_{\mathbf{g}}^{-1},G_1,\ldots,G_j)>0,
\]
where $P_j$ is a universal polynomial with nonnegative coefficients, such that
\begin{equation}
\label{eq:amann-buc-metrica-inversa}
 \|(\mathbf{g}_\kappa^{-1})^{ab}\|_{BUC^j(Q_\kappa)}\leq G'_j
\end{equation}
for all charts and indices. The same transformation rule, the Faà di Bruno formula of Theorem~\ref{faa di bruno multivariable}, and the bounds on mixed transitions show that the definition is independent of the representative atlas.

\begin{theorem}[From bounded geometry to uniform regularity]
\label{teo:geometria-acotada-implica-regularidad-uniforme-amann}
Every Riemannian manifold without boundary and with bounded geometry admits a uniformly regular structure. The atlas constants can be chosen to depend on $m$, a common normal radius, and the curvature bounds in Definition~\ref{def: geometria acotada}, never on the chart.
\end{theorem}

\begin{proof}
Let $i_0:=\operatorname{inj}(M)>0$, and let $r_0>0$ be the radius in Theorem~\ref{teo: equivalencias de geometria acotada}. Choose $a>0$ so that
\[
 4a\sqrt m<\min\{i_0,r_0\},
 \qquad \delta:=\frac a4,
 \qquad \rho:=\frac12.
\]
By the maximality argument of Theorem~\ref{thm:trivializacion-geodesica}, there exists a family $(p_i)_{i\in I}$ that is $\delta$-separated and whose balls of radius $\delta$ cover $M$. Choose isometries $B_i\colon \mathbb R^m\longrightarrow T_{p_i}M$ and define
\[
 \kappa_i^{-1}(x):=\exp_{p_i}(B_i(ax)),
 \qquad x\in Q^m.
\]
The charts are well defined because $aQ^m\subset B(0,a\sqrt m)$. If $p\in B_{\mathbf{g}}(p_i,\delta)$, then $p=\exp_{p_i}(B_iv)$ with $|v|<\frac{a}{4}$ and $\frac{v}{a}\in\rho Q^m$. Thus, the shrunken charts cover $M$.

If $U_i\cap U_j\neq\varnothing$, then $d_{\mathbf{g}}(p_i,p_j)<2a\sqrt m$. The balls $B_{\mathbf{g}}(p_j,\frac{\delta}{2})$ are pairwise disjoint and, for every index meeting $U_i$, lie in $B_{\mathbf{g}}(p_i,2a\sqrt m+\frac{\delta}{2})$. With the constants from Lemma~\ref{lema: volumen bolas geometria acotada}, we obtain
\begin{equation}
\label{eq:multiplicidad-atlas-amann-explicita}
 \operatorname{mult}(\mathfrak K)
 \leq N_0
 :=\left\lceil
 \frac{v_+(2a\sqrt m+\frac{\delta}{2})}{v_-(\frac{\delta}{2})}
 \right\rceil.
\end{equation}

Let $\phi_i=B_i^{-1}\circ\exp_{p_i}^{-1}$ be the unscaled normal chart. If $U_i\cap U_j\neq\varnothing$, then $d_{\mathbf{g}}(p_i,p_j)<2a\sqrt m$ and, for every $x\in Q^m$,
\[
 d_{\mathbf{g}}\bigl(p_j,\exp_{p_i}(B_i(ax))\bigr)<3a\sqrt m.
\]
The choice $\displaystyle 4a\sqrt m<\min\{i_0,r_0\}$ therefore allows the coordinate transition to be extended to all of $Q^m$ by
\[
 T_{ji}^{\mathrm{ext}}(x)
 :=a^{-1}B_j^{-1}\exp_{p_j}^{-1}
       \bigl(\exp_{p_i}(B_i(ax))\bigr).
\]
Theorem~\ref{teo: equivalencias de geometria acotada} provides constants $A_\ell>0$, independent of $i,j$, for these exponential comparisons at the preceding radii. Since, on the actual domain of the transition,
\[
 \kappa_j\circ\kappa_i^{-1}(x)
 =a^{-1}(\phi_j\circ\phi_i^{-1})(ax),
\]
we have, for $\ell\geq1$,
\begin{equation}
\label{eq:cota-transiciones-amann-reescaladas}
 \displaystyle\sup_{x\in\kappa_i(U_i\cap U_j)}
 \|D^\ell(\kappa_j\circ\kappa_i^{-1})(x)\|_{\operatorname{op}}
 \leq a^{\ell-1}A_\ell,
\end{equation}
and the order-zero bound for the actual transition is $\sqrt m$ because its image lies in $Q^m$. Since $Q^m$ is convex, the estimate of order $\ell+1$ for $T_{ji}^{\mathrm{ext}}$ and the mean value theorem give
\[
 \|D^\ell T_{ji}^{\mathrm{ext}}(x)
       -D^\ell T_{ji}^{\mathrm{ext}}(y)\|_{\operatorname{op}}
 \leq a^\ell A_{\ell+1}|x-y|,
 \qquad x,y\in Q^m.
\]
Restricting to each overlap domain gives the required uniform continuity; the constants depend only on $a,m,A_1,\ldots,A_{\ell+1}$.

Denote by $C_\ell$ the bounds for the metric coefficients and their inverses in normal charts of radius $2a\sqrt m$. The transformation rule under the rescaling $x\mapsto ax$ is, componentwise,
\[
 (g_{\kappa_i})_{bc}(x)
 =a^2(g_{\phi_i})_{bc}(ax).
\]
The bounds on the entries of $(g_{\phi_i})_{bc}$ and $(g_{\phi_i})^{bc}$ imply, by the spectral argument of Lemma~\ref{lema: equivalencia uniforme metricas}, that the eigenvalues of $(g_{\phi_i})_{bc}$ lie in $[(mC_0)^{-1},mC_0]$. Consequently,
\[
 c_{\mathbf{g}}:=\frac{a^2}{mC_0},
 \qquad C_{\mathbf{g}}:=a^2mC_0
\]
satisfy \eqref{eq:amann-comparacion-metrica-explicita}. Moreover,
\[
 \|D^\ell(\mathbf{g}_{\kappa_i})_{ab}\|_{L^\infty(Q^m)}
 \leq a^{\ell+2}C_\ell,
 \qquad 0\leq\ell\leq j.
\]
Summing these bounds provides $\displaystyle G_j:=\displaystyle\sum_{\ell=0}^ja^{\ell+2}C_\ell$ in \eqref{eq:amann-buc-metrica}; the bound of order $j+1$ gives uniform continuity of the derivative of order $j$. This verifies all axioms for the uniformly regular atlas and metric.
\end{proof}

\begin{theorem}[Full comparison in the absence of boundary]
\label{teo:amann-ur-equivale-geometria-acotada}
A Riemannian manifold without boundary is uniformly regular if and only if it has bounded geometry in the sense of Definition~\ref{def: geometria acotada}.
\end{theorem}

\begin{proof}
The implication from bounded geometry is the preceding theorem. For the converse, let $\mathfrak K$ be a uniformly regular atlas with shrinking factor $\rho\in(0,1)$ and the metric bounds in \eqref{eq:amann-comparacion-metrica-explicita}. Work in the normalized charts $Q^m=(-1,1)^m$; there are no boundary charts here. The formula
\[
 \Gamma^a_{bc}
 =\frac12\displaystyle\sum_{d=1}^m g^{ad}
   (\partial_b g_{cd}+\partial_c g_{bd}-\partial_d g_{bc})
\]
and \eqref{eq:recurrencia-derivadas-metrica-inversa-amann} provide constants $A_0,A_1\geq0$ such that
\[
 \|\Gamma(x)\|_{\mathrm{op}}\leq A_0,
 \qquad \|D\Gamma(x)\|_{\mathrm{op}}\leq A_1,
 \qquad x\in Q^m,
\]
uniformly over the charts. Regard $\Gamma(x)$ as a bilinear form with values in $\mathbb R^m$ and use Euclidean norms.

Set $\delta=1-\rho>0$. Choose $a>0$ independent of the chart and small enough that
\[
 2a<\frac\delta2,
 \qquad 4A_0a\leq\frac12,
 \qquad b:=4A_1a^2+4A_0a\leq\frac14.
\]
Let $p\in M$. The shrinking property allows us to choose a chart $\kappa$ with $x=\kappa(p)\in\rho Q^m$. For $v\in\mathbb R^m$, $|v|<a$, solve the coordinate geodesic equation
\[
 z''(t)=-\Gamma(z(t))(z'(t),z'(t)),
 \qquad z(0)=x,\quad z'(0)=v.
\]
As long as $0\leq t\leq1$ and $|z'(t)|\leq2a$, its integral form gives
\[
 |z'(t)|\leq a+4A_0a^2\leq\frac32a,
 \qquad |z(t)-x|\leq2a<\frac\delta2.
\]
By continuity, the speed cannot reach $2a$ before $t=1$. The position remains at distance at least $\delta/2$ from the boundary of the cube, and the velocity remains bounded. Thus, the solution stays in a compact subset of the domain of the first-order equation for $(z,z')$, so Theorem~\ref{teo:criterio-prolongacion-edo-banach} extends it to all of $[0,1]$. In particular, the exponential map is defined on these initial vectors.

Fix $w\in\mathbb R^m$ and denote by $Y(t)=D_vz(t)[w]$ the variation with respect to the initial velocity. Smooth dependence of ODE solutions on their data, Theorem~\ref{teo:dependencia-parametros-edo-banach}, and differentiation of the equation give
\[
 \begin{split}
 Y''(t)={}&-D\Gamma(z(t))[Y(t)](z'(t),z'(t))
          -2\Gamma(z(t))(Y'(t),z'(t)),\\
 Y(0)={}&0,\qquad Y'(0)=w.
 \end{split}
\]
Let $\displaystyle S=\displaystyle\max_{0\leq t\leq1}|Y'(t)|$. Since $|Y(t)|\leq S$, the integrated equation gives
\[
 |Y'(t)-w|\leq(4A_1a^2+4A_0a)S=bS.
\]
Thus, $S\leq|w|+bS$, whence
\[
 S\leq\frac{|w|}{1-b},
 \qquad |Y(1)-w|\leq\frac{b}{1-b}|w|\leq\frac13|w|.
\]
If $F_x(v)=z(1)$, we have obtained
$\|DF_x(v)-I\|_{\mathrm{op}}\leq1/3$
on $B(0,a)$. This domain is convex; integrating the derivative along the segment between $v$ and $w$ gives
\[
 |F_x(v)-F_x(w)-(v-w)|\leq\frac13|v-w|,
 \qquad |F_x(v)-F_x(w)|\geq\frac23|v-w|.
\]
Thus, $F_x$ is injective and its derivative is invertible. The inverse function theorem makes it a diffeomorphism onto its image. Since $F_x$ is the coordinate expression of $\exp_p$ and
\[
 |d\kappa_p(\mathbf v)|
 \leq c_{\mathbf g}^{-1/2}|\mathbf v|_{\mathbf g_p},
 \qquad \mathbf v\in T_pM,
\]
the exponential map is a diffeomorphism onto its image on the tangent ball of radius $a\sqrt{c_{\mathbf g}}$. This bound is independent of $p$ and proves
\[
 \operatorname{inj}(M)\geq a\sqrt{c_{\mathbf g}}>0.
\]
Completeness of the components now follows from Remark~\ref{obs: geometria acotada implica completa}.

The curvature bounds remain. In each chart,
\[
 R_{ijb}{}^a
 =\partial_i\Gamma^a_{jb}-\partial_j\Gamma^a_{ib}
  +\displaystyle\sum_{d=1}^m
   (\Gamma^d_{jb}\Gamma^a_{id}-\Gamma^d_{ib}\Gamma^a_{jd}).
\]
Every coordinate derivative of order $k$ of these components is a finite sum of products of derivatives of $\Gamma$ of orders at most $k+1$. These are uniformly bounded by the metric bounds through order $k+2$ and the recurrence for its inverse. To pass to covariant derivatives, use the tensor formula: add to the coordinate derivative a term $\Gamma T$ for each contravariant index and subtract one for each covariant index. At order zero, we have just bounded $R$. If the components of $\nabla^jR$ have been expressed as sums of products of coordinate derivatives of $R$ and $\Gamma$, applying this formula and the Leibniz rule produces an expression of the same type for $\nabla^{j+1}R$, with orders increased by at most one. Induction provides uniform bounds for all these components. Finally, \eqref{eq:amann-comparacion-metrica-explicita} and the inverse-metric bounds convert component bounds into the tensor bounds
$\displaystyle\sup_{p\in M}|\nabla^k\mathbf R(p)|_{\mathbf g}<\infty$
for every $k\in\mathbb N_0$. Both bounded-geometry conditions hold.
\end{proof}

The definitions and this equivalence correspond to \cite[Definitions~3.1--3.3 and Theorem~3.5]{Amann2025FunctionSpaces}.

\begin{proposition}[Amann's localization system]
\label{prop:sistema-localizacion-amann}
\index{Amann's localization system}
Let $M$ be a smooth manifold with or without boundary, and let $\mathfrak K$ be a uniformly regular atlas of $M$ with shrinking factor $\rho$. There exist functions $\pi_\kappa\in C_c^\infty(U_\kappa,[0,1])$ and a function $\chi\in C_c^\infty(Q^m,[0,1])$ such that
\[
 \displaystyle\sum_{\kappa\in K}\pi_\kappa^2=1,
 \qquad
 \chi=1\quad\hbox{on }
 \operatorname{supp}(\pi_\kappa\circ\kappa^{-1}),
\]
and, for each $j\in\mathbb N_0$, there exists $D_j>0$, independent of $\kappa$, with
\[
 \|\pi_\kappa\circ\kappa^{-1}\|_{BUC^j(Q_\kappa)}\leq D_j.
\]
If $\mathfrak K$ is the atlas constructed in Theorem~\ref{teo:geometria-acotada-implica-regularidad-uniforme-amann}, or a uniformly regular atlas equivalent to it, this is also an admissible quadratic localization system in the sense of Definition~\ref{def:sistema-cuadratico-localizacion-admisible}.
\end{proposition}

\begin{proof}
Choose $\rho<\theta<\theta'<1$. For $0<a<b<1$, let
\[
 \eta_{a,b}(t):=
 \frac{\vartheta(b^2-t^2)}
 {\vartheta(b^2-t^2)+\vartheta(t^2-a^2)},
 \qquad t\in\mathbb R,
\]
where $\vartheta$ is the flat function used in the proof of Theorem~\ref{particionesdelaunidad}. If $b_0=(\rho+\theta)/2$, define
\[
 \psi(x_1,\dots,x_m):=\prod_{a=1}^m\eta_{\rho,b_0}(x_a).
\]
Then $\psi\in C_c^\infty(Q^m,[0,1])$, it equals one on $\rho Q^m$, and its support lies in $[-b_0,b_0]^m\subset\theta Q^m$. On a half-box, use its restriction. Let $\psi_\kappa=\psi\circ\kappa$ on $U_\kappa$, extended by zero, and set
\[
 S(p):=\left(\displaystyle\sum_{\kappa\in K}\psi_\kappa(p)^2\right)^{\frac{1}{2}},
 \qquad \pi_\kappa:=\frac{\psi_\kappa}{S}.
\]
The family of supports entering this sum is not merely locally finite: the number of these supports meeting any given chart has a uniform bound. We prove both assertions at once. Set
\[
 N:=\operatorname{mult}(\mathfrak K),
 \qquad
 \varepsilon
 :=\frac{1-\theta}{2\max\{1,C_{\mathrm{tr}}(1)\}}>0.
\]
Fix a chart $\kappa_0$ and consider the indices $\kappa$ for which $\operatorname{supp}\psi_\kappa\cap U_{\kappa_0}\neq\varnothing$. For each one, choose $q_\kappa\in\operatorname{supp}\psi_\kappa\cap U_{\kappa_0}$. Then $\kappa(q_\kappa)\in\operatorname{supp}\psi\subset\theta Q_\kappa$, whose distance from the nongeometric faces of $Q_\kappa$ is at least $1-\theta$. The bound $BUC^1$ for $\kappa\circ\kappa_0^{-1}$ and the first-exit-time argument, applied to each segment starting at $\kappa_0(q_\kappa)$, show that
\begin{equation}
\label{eq:bola-uniforme-traslape-cortes-amann}
 B\bigl(\kappa_0(q_\kappa),\varepsilon\bigr)\cap Q_{\kappa_0}
 \subseteq \kappa_0(U_\kappa\cap U_{\kappa_0}).
\end{equation}
On a half-box, a path contained in $Q_{\kappa_0}$ cannot leave through the geometric face: chart transitions preserve the boundary. Thus, only distance from the nongeometric faces enters the argument.

All sets on the left-hand side of \eqref{eq:bola-uniforme-traslape-cortes-amann} lie in the bounded set
\[
 W:=\{x\in\mathbb R^m\mid \operatorname{dist}(x,Q_{\kappa_0})<\varepsilon\}.
\]
Since $Q_{\kappa_0}$ is a box or a half-box, there exists $c_m>0$, depending only on the dimension, such that
\[
 \lambda_m\bigl(B(z,\varepsilon)\cap Q_{\kappa_0}\bigr)
 \geq c_m\varepsilon^m,
 \qquad z\in Q_{\kappa_0}.
\]
On the other hand, each point of $U_{\kappa_0}$ belongs to at most $N$ charts. Summing the indicator functions of the preceding balls and integrating over $W$ gives
\[
 \#\{\kappa\mid
 \operatorname{supp}\psi_\kappa\cap U_{\kappa_0}\neq\varnothing\}
 c_m\varepsilon^m
 \leq N\lambda_m(W).
\]
The measure of $W$ has a bound depending only on $m$ and $\varepsilon$, not on $\kappa_0$. Thus, there exists an integer $N_1$, independent of the chart, bounding the cardinality of the preceding index set. This proves, in particular, local finiteness of the family of supports.

Since the shrunken charts cover, $1\leq S\leq\sqrt N$. In chart $\kappa$, at most $N_1$ functions of the form $\psi\circ\widetilde\kappa\circ\kappa^{-1}$ contribute. The Faà di Bruno formula of Theorem~\ref{faa di bruno multivariable}, the transition bounds, and the Leibniz rule bound all their $BUC^j$ norms by constants independent of the indices. If $\displaystyle F_\kappa:=\displaystyle\sum_{\widetilde\kappa\in K}\psi_{\widetilde\kappa}^2$ is expressed in chart $\kappa$, then $1\leq F_\kappa\leq N$ and $S^{-1}=F_\kappa^{-1/2}$. For $|\alpha|\geq1$, the same formula is a finite sum of terms of the form
\[
 c_{\alpha,q,\beta}\,
 F_\kappa^{-\frac12-q}
 \prod_{\nu=1}^qD^{\beta_\nu}F_\kappa,
 \qquad
 \beta_\nu\neq0,\quad
 \beta_1+\cdots+\beta_q=\alpha,
\]
with $1\leq q\leq|\alpha|$. The lower bound for $F_\kappa$ and its derivative bounds explicitly control every product, uniformly over the charts. Applying the Leibniz rule to $\pi_\kappa=\psi_\kappa S^{-1}$ gives the required bounds. Moreover, $\displaystyle \sum_{\kappa\in K}\pi_\kappa^2=S^{-2}\displaystyle\sum_{\kappa\in K}\psi_\kappa^2=1$.

If $b_1=(\theta+\theta')/2$, take
\[
 \chi(x_1,\dots,x_m):=\prod_{a=1}^m\eta_{\theta,b_1}(x_a).
\]
Then $\chi\in C_c^\infty(Q^m,[0,1])$, it equals one on $\theta Q^m$, and its support lies in $\theta'Q^m$. Therefore, $\chi=1$ on all coordinate supports. When $\mathfrak K$ is the rescaled normal atlas of the preceding theorem, the transition bounds already established verify \ref{item:B1}; the bounds just obtained for the cutoffs verify the remaining conditions of Definition~\ref{def:sistema-cuadratico-localizacion-admisible}. For an equivalent atlas, the bounds on mixed transitions transfer these conditions through the same Faà di Bruno formula and Leibniz rule.
\end{proof}

The preceding system gives the following localization operators. Taking $h_i=\pi_{\kappa_i}$ and $\chi_i:=\chi\circ\kappa_i$, the elementary operators on smooth functions are
\[
 \mathcal S_0u=
 \left(\widetilde{(h_i u)\circ\phi_i^{-1}}\right)_{i\in I},
 \qquad
 (\mathcal R_0v)(p)=
 \displaystyle\sum_{\{i\in I\mid p\in U_i\}}
 h_i(p)\chi_i(p)v_i(\phi_i(p)).
\]
The second sum is locally finite. Moreover, $\chi_i=1$ on $\operatorname{supp}h_i$ and $\displaystyle \sum_{i\in I} h_i^2=1$; hence, $\mathcal R_0\mathcal S_0u=u$. These formulas are used here only on smooth functions. Distributional reconstruction and the Sobolev norm bounds are established after these spaces are defined, using the admissibility conditions verified for this system.

The differential part of \ref{item:B1} admits the following metric criterion. The overlap hypothesis is kept separate because it does not follow from estimates on the metric coefficients.
\begin{lemma}[Metric criterion for admissibility]\label{lem:geometria-acotada-variedad-riemanniana-trivializacion-uniformemente-localmente-finit}
Let $(M,\mathbf{g})$ be a Riemannian manifold without boundary with $\operatorname{inj}(M)>0$, and let \(\mathcal T=(U_i,\phi_i,h_i)_{i\in I}\) be a uniformly locally finite trivialization. Assume also that the chart images are uniformly normalized: there exists $R>0$ such that $\phi_i(U_i)\subseteq B_{\mathrm{euc}}(0,R)$ for every $i\in I$. Suppose that its cover satisfies the two cross-overlap bounds in \ref{item:B1} relative to a uniform geodesic atlas with the separation and covering properties specified there, chosen with sufficiently small radius. Let \(g_{ab}^{(i)}\) and \(g_{(i)}^{ab}\) be the coefficients of \(\mathbf{g}\) and its inverse in coordinates \(\phi_i\). Then \((M,\mathbf{g})\) has bounded geometry and \(\mathcal T\) satisfies \ref{item:B1} if and only if, for each \(k\in\mathbb N_0\), there exists \(C_k>0\) such that, for \(|\alpha|\leq k\),
\[
 |D^{\alpha}g_{ab}^{(i)}|\leq C_{k}\qquad \text{and}\qquad |D^{\alpha}g^{ab}_{(i)}|\leq C_{k}
 \]
in every chart \((U_i,\phi_i)\).
\end{lemma}

\begin{proof}
\noindent\textit{($\Leftarrow$)} Suppose that, for every $k\in \mathbb{N}_{0}$, there exists $C_{k}>0$ such that, for every chart $(U_i,\phi_i)$ of the trivialization $\mathcal{T}$,
\[
|D^{\alpha}g_{ab}^{(i)}|\leq C_{k}
\qquad\text{and}\qquad
|D^{\alpha}g_{(i)}^{ab}|\leq C_{k},
\qquad \forall |\alpha|\leq k.
\]

First, we prove that $(M,\mathbf{g})$ has bounded geometry. Since $\operatorname{inj}(M)>0$ by hypothesis, it suffices to show that, for each $m\in \mathbb{N}_{0}$, there exists a constant $C_m'>0$ such that
\[
|\nabla^{m}R|_{\mathbf{g}}\leq C_m'
\]
throughout $M$.

Fix a chart $(U_i,\phi_i)$ of $\mathcal{T}$ and write $\phi_i=(x^1,\dots,x^n)$.

To simplify notation, denote simply by \[
g_{ab},\qquad g^{ab},\qquad \Gamma_{ab}^{c},\qquad R_{abc}{}^{d}
\] the components of the metric, inverse metric, Christoffel symbols, and curvature tensor in these coordinates, which of course depend on $i$.

By Corollary~\ref{expresion local christoffel levi civita}, the local expression for the Christoffel symbols of the Levi--Civita connection is:
\[
\Gamma_{ab}^{c}
=
\frac{1}{2}g^{c\ell}
\bigl(
\partial_{a}g_{b\ell}
+
\partial_{b}g_{a\ell}
-
\partial_{\ell}g_{ab}
\bigr).
\]
Applying the Leibniz rule and the hypothesis on the partial derivatives of $g_{ab}$ and $g^{ab}$, we conclude that, for each $k\in\mathbb{N}_{0}$, there exists a constant $C_k^{(1)}>0$, independent of the chart, such that
\[
|D^\alpha \Gamma_{ab}^{c}|\le C_k^{(1)},
\qquad |\alpha|\le k.
\]

Now using the local expression for the curvature tensor (Proposition~\ref{1,3 curvatura en coordenadas}),
\[
R_{abc}{}^{\ell}
=
\partial_a \Gamma_{bc}^{\ell}
-
\partial_b \Gamma_{ac}^{\ell}
+
\Gamma_{bc}^{m}\Gamma_{am}^{\ell}
-
\Gamma_{ac}^{m}\Gamma_{bm}^{\ell},
\]
the Leibniz rule again gives, for each $k\in\mathbb{N}_{0}$, a constant $C_k^{(2)}>0$, independent of the chart, such that
\[
|D^\alpha R_{abc}{}^{\ell}|\le C_k^{(2)},
\qquad |\alpha|\le k.
\]

Finally, by Lemma~\ref{lem:local-expression-higher-order}, the coordinate components of $\nabla^{m}R$ are finite linear combinations of partial derivatives of the components of $R$ of order at most $m$, with coefficients polynomial in the Christoffel symbols and their partial derivatives of order at most $m-1$. Since all these quantities are uniformly bounded, for each $m\in \mathbb{N}_{0}$ there exists a constant $C_m''>0$ such that the components of $\nabla^mR$ are uniformly bounded in every chart of $\mathcal{T}$.

Moreover, the uniform bounds on the components of the metric and its inverse imply that, for the fixed tensor type of $\nabla^mR$, there exists $C_{m,\mathbf{g}}>0$, independent of the chart, such that the pointwise norm induced by $\mathbf{g}$ does not exceed $C_{m,\mathbf{g}}$ times the Euclidean norm of its components. Consequently, for each $m\in\mathbb{N}_{0}$, there exists a constant $C_m'>0$ such that
\[
|\nabla^{m}R|_{\mathbf{g}}\leq C_m'
\]
throughout $M$. This proves that $(M,\mathbf{g})$ has bounded geometry.

We now prove that $\mathcal{T}$ satisfies \ref{item:B1}. Let \(r_0>0\) be a radius for which condition (2) of Theorem~\ref{teo: equivalencias de geometria acotada} holds, and let \[
\mathcal{A}^{\operatorname{geo}}
=
(U_j^{\operatorname{geo}},\phi_j^{\operatorname{geo}})_{j\in J}
\] be the geodesic atlas in the hypothesis. Here ``sufficiently small radius'' means that its common radius satisfies \(0<r\leq \frac{r_0}{10}\). Fix $i\in I$ and $j\in J$ such that
\[
U_i\cap U_j^{\operatorname{geo}}\neq \varnothing.
\]
Write
\[
\phi_i=(x^1,\dots,x^n),
\qquad
\phi_j^{\operatorname{geo}}=(y^1,\dots,y^n).
\]
Denote by \[
\widetilde g_{ab},\qquad \widetilde g^{ab},\qquad \widetilde\Gamma_{ab}^{c}
\] the components of the metric and its inverse and the Christoffel symbols in geodesic coordinates \(y\).

Since we have already proved that $(M,\mathbf{g})$ has bounded geometry, Theorem~\ref{teo: equivalencias de geometria acotada} implies that, in the uniform geodesic atlas, the components \(\widetilde g_{ab}\), \(\widetilde g^{ab}\), and all their partial derivatives are uniformly bounded. In particular, for each $k\in\mathbb{N}_{0}$, there exists $A_k>0$ such that
\[
|D^\alpha \widetilde g_{ab}|\le A_k
\qquad\text{and}\qquad
|D^\alpha \widetilde g^{ab}|\le A_k,
\qquad \forall |\alpha|\le k.
\]

Order zero is controlled by normalization: the image of $\phi_i\circ(\phi_j^{\operatorname{geo}})^{-1}$ lies in $B_{\mathrm{euc}}(0,R)$ and that of its inverse in $B_{\mathrm{euc}}(0,r)$. Metric bounds do not replace this hypothesis, since translating the origin of a chart leaves its metric coefficients unchanged. We now begin with first derivatives of the coordinate transition
\[
\phi_i\circ (\phi_j^{\operatorname{geo}})^{-1}
\qquad\text{and}\qquad
\phi_j^{\operatorname{geo}}\circ \phi_i^{-1}.
\]

By hypothesis, the components of the metric and its inverse in chart $(U_i,\phi_i)$ are uniformly bounded. Since $(M,\mathbf{g})$ also has bounded geometry, the same holds in the geodesic chart $(U_j^{\operatorname{geo}},\phi_j^{\operatorname{geo}})$. By Lemma~\ref{lema: equivalencia uniforme metricas}, there exist positive constants $c_0,C_0$, independent of $i$, $j$, and the point considered, such that
\[
c_0\|d\phi_i(v)\|^2
\le
|v|_{\mathbf{g}}^2
\le
C_0\|d\phi_i(v)\|^2
\]
for every \(v\in TM\restriction_{U_i}\), and
\[
c_0\|d\phi_j^{\operatorname{geo}}(v)\|^2
\le
|v|_{\mathbf{g}}^2
\le
C_0\|d\phi_j^{\operatorname{geo}}(v)\|^2
\]
for every \(v\in TM\restriction_{U_j^{\operatorname{geo}}}\).

Now let
\[
F_{ij}:=\phi_i\circ (\phi_j^{\operatorname{geo}})^{-1}.
\]
Fix \(p\in U_i\cap U_j^{\operatorname{geo}}\) and let \(\xi\in\mathbb{R}^{n}\). Define
\[
v:=(d\phi_j^{\operatorname{geo}}|_{p})^{-1}(\xi)\in T_pM.
\]
Then
\[
d\phi_i|_{p}(v)
=
d(\phi_i\circ(\phi_j^{\operatorname{geo}})^{-1})|_{\phi_j^{\operatorname{geo}}(p)}(\xi)
=
DF_{ij}\bigl(\phi_j^{\operatorname{geo}}(p)\bigr)\,\xi.
\]
Applying the preceding inequalities to this vector \(v\) gives
\[
c_0\bigl\|DF_{ij}\bigl(\phi_j^{\operatorname{geo}}(p)\bigr)\xi\bigr\|^2
\le
|v|_{\mathbf{g}}^2
\le
C_0\|\xi\|^2.
\]
Therefore,
\[
\bigl\|DF_{ij}\bigl(\phi_j^{\operatorname{geo}}(p)\bigr)\xi\bigr\|
\le
\sqrt{\frac{C_0}{c_0}}\|\xi\|,
\]
and hence
\[
\bigl\|DF_{ij}\bigl(\phi_j^{\operatorname{geo}}(p)\bigr)\bigr\|_{\operatorname{op}}
\le
\sqrt{\frac{C_0}{c_0}}.
\]
Since the operator norm of $DF_{ij}$ is uniformly bounded, the partial derivatives of order $1$ are uniformly bounded.

For the inverse transition, define
\[
G_{ij}:=\phi_j^{\operatorname{geo}}\circ \phi_i^{-1},
\]
Now taking \(v=(d\phi_i|_p)^{-1}(\xi)\) and applying the preceding two metric inequalities with the charts interchanged gives, for every \(p\in U_i\cap U_j^{\operatorname{geo}}\),
\[
\bigl\|DG_{ij}\bigl(\phi_i(p)\bigr)\bigr\|_{\operatorname{op}}
\le
\sqrt{\frac{C_0}{c_0}}.
\]

We turn to second derivatives. Applying item (c) of Lemma~\ref{lema: transformacion metrica Christoffel} to the charts \((U_i,\phi_i)\) and \((U_j^{\operatorname{geo}},\phi_j^{\operatorname{geo}})\) gives, for every \(p\in U_i\cap U_j^{\operatorname{geo}}\),
\[
\Gamma_{rs}^{k}(p)
=
\frac{\partial x^{k}}{\partial y^{c}}(p)
\frac{\partial y^{a}}{\partial x^{r}}(p)
\frac{\partial y^{b}}{\partial x^{s}}(p)\,
\widetilde\Gamma_{ab}^{c}(p)
+
\frac{\partial x^{k}}{\partial y^{c}}(p)
\frac{\partial^{2}y^{c}}{\partial x^{r}\partial x^{s}}(p).
\]
Multiplying by \(\displaystyle\frac{\partial y^{m}}{\partial x^{k}}(p)\) and using the chain rule, we conclude that
\[
\frac{\partial^{2}y^{m}}{\partial x^{r}\partial x^{s}}(p)
=
\frac{\partial y^{m}}{\partial x^{k}}(p)\,\Gamma_{rs}^{k}(p)
-
\frac{\partial y^{a}}{\partial x^{r}}(p)
\frac{\partial y^{b}}{\partial x^{s}}(p)\,
\widetilde\Gamma_{ab}^{m}(p).
\]
Since \[
y^m
=
\pi_m\circ \phi_j^{\operatorname{geo}}
=
\pi_m\circ G_{ij}\circ \phi_i
\] on \(U_i\cap U_j^{\operatorname{geo}}\), we obtain
\[
\frac{\partial^{2}G_{ij}^{m}}{\partial x^{r}\partial x^{s}}\bigl(\phi_i(p)\bigr)
=
\frac{\partial G_{ij}^{m}}{\partial x^{k}}\bigl(\phi_i(p)\bigr)\,\Gamma_{rs}^{k}(p)
-
\frac{\partial G_{ij}^{a}}{\partial x^{r}}\bigl(\phi_i(p)\bigr)
\frac{\partial G_{ij}^{b}}{\partial x^{s}}\bigl(\phi_i(p)\bigr)\,
\widetilde\Gamma_{ab}^{m}(p).
\]
The first derivatives of \(G_{ij}\) are already known to be uniformly bounded, as are the components of \(\Gamma\) and \(\widetilde\Gamma\). Thus, the second derivatives of \(G_{ij}\) are uniformly bounded as well.

Interchanging the roles of the two charts in item (c) of Lemma~\ref{lema: transformacion metrica Christoffel} gives, for every \(p\in U_i\cap U_j^{\operatorname{geo}}\),
\[
\widetilde\Gamma_{ab}^{c}(p)
=
\frac{\partial y^{c}}{\partial x^{k}}(p)
\frac{\partial x^{r}}{\partial y^{a}}(p)
\frac{\partial x^{s}}{\partial y^{b}}(p)\,
\Gamma_{rs}^{k}(p)
+
\frac{\partial y^{c}}{\partial x^{k}}(p)
\frac{\partial^{2}x^{k}}{\partial y^{a}\partial y^{b}}(p).
\]
Multiplying by \(\displaystyle\frac{\partial x^{m}}{\partial y^{c}}(p)\) gives
\[
\frac{\partial^{2}x^{m}}{\partial y^{a}\partial y^{b}}(p)
=
\frac{\partial x^{m}}{\partial y^{c}}(p)\,\widetilde\Gamma_{ab}^{c}(p)
-
\frac{\partial x^{r}}{\partial y^{a}}(p)
\frac{\partial x^{s}}{\partial y^{b}}(p)\,
\Gamma_{rs}^{m}(p).
\]
Since \[
x^m
=
\pi_m\circ \phi_i
=
\pi_m\circ F_{ij}\circ \phi_j^{\operatorname{geo}}
\] on \(U_i\cap U_j^{\operatorname{geo}}\), we obtain
\[
\frac{\partial^{2}F_{ij}^{m}}{\partial y^{a}\partial y^{b}}\bigl(\phi_j^{\operatorname{geo}}(p)\bigr)
=
\frac{\partial F_{ij}^{m}}{\partial y^{c}}\bigl(\phi_j^{\operatorname{geo}}(p)\bigr)\,\widetilde\Gamma_{ab}^{c}(p)
-
\frac{\partial F_{ij}^{r}}{\partial y^{a}}\bigl(\phi_j^{\operatorname{geo}}(p)\bigr)
\frac{\partial F_{ij}^{s}}{\partial y^{b}}\bigl(\phi_j^{\operatorname{geo}}(p)\bigr)\,
\Gamma_{rs}^{m}(p).
\]
Thus, the second derivatives of \(F_{ij}\) are also uniformly bounded.

We now proceed by induction on the order of differentiation. Suppose that, for some \(k\geq 2\), all partial derivatives of \(F_{ij}\) and \(G_{ij}\) of order at most \( k\) are already known to be uniformly bounded. We will show that those of order \(k+1\) are uniformly bounded as well.

Consider the identity
\[
\frac{\partial^{2}F_{ij}^{m}}{\partial y^{a}\partial y^{b}}\bigl(\phi_j^{\operatorname{geo}}(p)\bigr)
=
\frac{\partial F_{ij}^{m}}{\partial y^{c}}\bigl(\phi_j^{\operatorname{geo}}(p)\bigr)\,\widetilde\Gamma_{ab}^{c}(p)
-
\frac{\partial F_{ij}^{r}}{\partial y^{a}}\bigl(\phi_j^{\operatorname{geo}}(p)\bigr)
\frac{\partial F_{ij}^{s}}{\partial y^{b}}\bigl(\phi_j^{\operatorname{geo}}(p)\bigr)\,
\Gamma_{rs}^{m}(p).
\]
Composing with \((\phi_j^{\operatorname{geo}})^{-1}\) gives, on the corresponding Euclidean open subset,
\[
\frac{\partial^{2}F_{ij}^{m}}{\partial y^{a}\partial y^{b}}
=
\frac{\partial F_{ij}^{m}}{\partial y^{c}}\,
(\widetilde\Gamma_{ab}^{c}\circ (\phi_j^{\operatorname{geo}})^{-1})
-
\frac{\partial F_{ij}^{r}}{\partial y^{a}}
\frac{\partial F_{ij}^{s}}{\partial y^{b}}\,
(\Gamma_{rs}^{m}\circ \phi_i^{-1}\circ F_{ij}).
\]
Now let \(\alpha\) be a multi-index with \(|\alpha|\leq k-1\). Apply \(D^\alpha\) to this identity. By the Leibniz rule, the first term on the right-hand side is a finite sum of products of derivatives of \(F_{ij}\) of order at most \(k\) and derivatives of \(\widetilde\Gamma_{ab}^{c}\circ (\phi_j^{\operatorname{geo}})^{-1}\) of order at most \(k-1\), all uniformly bounded.

For the second term, use the Leibniz rule again and, when differentiating the composition
\[
(\Gamma_{rs}^{m}\circ \phi_i^{-1})\circ F_{ij},
\]
apply the Faà di Bruno formula of Theorem~\ref{faa di bruno multivariable} in the form stated there. If
\[
H:=\Gamma_{rs}^{m}\circ\phi_i^{-1},
\qquad F_{ij}=(F_{ij}^1,\dots,F_{ij}^n),
\]
then, for $1\leq|\beta|\leq k-1$,
\begin{equation}
\label{eq:faa-di-bruno-transiciones-cap13}
\begin{split}
D^\beta(H\circ F_{ij})(y)
={}&
\displaystyle\sum_{\substack{\delta\in\mathbb N_0^n\\
1\leq|\delta|\leq|\beta|}}
\ \displaystyle\sum_{\substack{\gamma_1,\ldots,\gamma_{|\delta|}\in\mathbb N_0^n\\
\gamma_1+\cdots+\gamma_{|\delta|}=\beta\\
|\gamma_\nu|>0}}
\ \displaystyle\sum_{\substack{
\ell_1,\dots,\ell_{|\delta|}\in\{1,\dots,n\}\\
\#\{\nu\mid \ell_\nu=a\}=\delta_a}}
c_{\beta,\delta,\gamma,\ell}
\bigl(D^\delta H\bigr)(F_{ij}(y))
\prod_{\nu=1}^{|\delta|}
D^{\gamma_\nu}F_{ij}^{\ell_\nu}(y).
\end{split}
\end{equation}
The sum is finite, and its coefficients depend only on the multi-indices. Each summand satisfies $|\delta|\leq|\beta|\leq k-1$ and $1\leq|\gamma_\nu|\leq|\beta|\leq k-1$. Derivatives of $H$ are bounded by the bounds already obtained for the Christoffel symbols in chart $i$, and derivatives of $F_{ij}$ are bounded by the induction hypothesis. Consequently, every derivative \[
D^\beta\bigl((\Gamma_{rs}^{m}\circ \phi_i^{-1})\circ F_{ij}\bigr),
\qquad |\beta|\le k-1,
\] is bounded by a constant independent of $i,j$. After applying the Leibniz rule to the three factors of the second term, no derivative of $F_{ij}$ on the right-hand side has order greater than $k$. The left-hand side is $D^\alpha\partial_{y^a}\partial_{y^b}F_{ij}^m$; varying $a,b$ and the multi-indices $\alpha$ with $|\alpha|=k-1$ gives every derivative of order $k+1$. We conclude that they are uniformly bounded.

For the inverse transition, begin with the identity
\[
\frac{\partial^{2}G_{ij}^{m}}{\partial x^{r}\partial x^{s}}\bigl(\phi_i(p)\bigr)
=
\frac{\partial G_{ij}^{m}}{\partial x^{k}}\bigl(\phi_i(p)\bigr)\,\Gamma_{rs}^{k}(p)
-
\frac{\partial G_{ij}^{a}}{\partial x^{r}}\bigl(\phi_i(p)\bigr)
\frac{\partial G_{ij}^{b}}{\partial x^{s}}\bigl(\phi_i(p)\bigr)\,
\widetilde\Gamma_{ab}^{m}(p),
\]
which, after composition with \(\phi_i^{-1}\), takes the form
\[
\frac{\partial^{2}G_{ij}^{m}}{\partial x^{r}\partial x^{s}}
=
\frac{\partial G_{ij}^{m}}{\partial x^{k}}\,
(\Gamma_{rs}^{k}\circ \phi_i^{-1})
-
\frac{\partial G_{ij}^{a}}{\partial x^{r}}
\frac{\partial G_{ij}^{b}}{\partial x^{s}}\,
(\widetilde\Gamma_{ab}^{m}\circ (\phi_j^{\operatorname{geo}})^{-1}\circ G_{ij}),
\]
and shows that all derivatives of order \(k+1\) of \(G_{ij}\) are uniformly bounded as well. Here we use exactly \eqref{eq:faa-di-bruno-transiciones-cap13} with $H=\widetilde\Gamma_{ab}^{m}\circ
(\phi_j^{\operatorname{geo}})^{-1}$ and with $F_{ij}$ replaced by $G_{ij}$; the same bounds on $|\delta|$ and $|\gamma_\nu|$ ensure that only derivatives covered by the induction hypothesis occur.

By induction, for each \(k\in \mathbb{N}_{0}\) there exists a constant \(C_k>0\), independent of \(i\) and \(j\), such that
\[
|D^\alpha F_{ij}|\le C_k
\qquad\text{and}\qquad
|D^\alpha G_{ij}|\le C_k
\]
for every \(|\alpha|\le k\). Together with the two cross-overlap bounds assumed in the statement, this proves \ref{item:B1}.

\medskip

\noindent\textit{($\Rightarrow$)} Now suppose that $(M,\mathbf{g})$ has bounded geometry and condition \ref{item:B1} holds.

Let \[
\mathcal{A}^{\operatorname{geo}}
=
(U_j^{\operatorname{geo}},\phi_j^{\operatorname{geo}})_{j\in J}
\] be one of the uniform geodesic atlases for which \ref{item:B1} holds. Fix \(i\in I\) and \(j\in J\) such that
\[
U_i\cap U_j^{\operatorname{geo}}\neq \varnothing.
\]
Write
\[
\phi_i=(x^1,\dots,x^n),
\qquad
\phi_j^{\operatorname{geo}}=(y^1,\dots,y^n),
\]
and let
\[
F_{ij}:=\phi_i\circ (\phi_j^{\operatorname{geo}})^{-1},
\qquad
G_{ij}:=\phi_j^{\operatorname{geo}}\circ \phi_i^{-1}.
\]
By \ref{item:B1}, for each \(k\in\mathbb{N}_{0}\) there exists a constant \(C_k>0\) such that all partial derivatives of \(F_{ij}\) and \(G_{ij}\) of order at most \(k\) are uniformly bounded.

Since \((M,\mathbf{g})\) has bounded geometry, Theorem~\ref{teo: equivalencias de geometria acotada} implies that, in the uniform geodesic atlas, the components of the metric and its inverse have uniformly bounded derivatives of every order. Denote by \[
\widetilde g_{ab},
\qquad
\widetilde g^{ab}
\] the components of the metric and inverse metric in coordinates \(y\). Then, for each \(k\in\mathbb{N}_{0}\), there exists \(A_k>0\) such that
\[
|D^\alpha \widetilde g_{ab}|\le A_k
\qquad\text{and}\qquad
|D^\alpha \widetilde g^{ab}|\le A_k
\]
for every \(|\alpha|\le k\).

Applying item (a) of Lemma~\ref{lema: transformacion metrica Christoffel} to the charts \((U_i,\phi_i)\) and \((U_j^{\operatorname{geo}},\phi_j^{\operatorname{geo}})\) gives, for every \(p\in U_i\cap U_j^{\operatorname{geo}}\),
\[
g_{rs}(p)
=
\frac{\partial y^{a}}{\partial x^{r}}(p)
\frac{\partial y^{b}}{\partial x^{s}}(p)\,
\widetilde g_{ab}(p).
\]
Since \[
y^a
=
\pi_a\circ G_{ij}\circ \phi_i
\] on \(U_i\cap U_j^{\operatorname{geo}}\), composing with \(\phi_i^{-1}\) gives
\[
g_{rs}\circ \phi_i^{-1}
=
\frac{\partial G_{ij}^{a}}{\partial x^{r}}
\frac{\partial G_{ij}^{b}}{\partial x^{s}}\,
\bigl(\widetilde g_{ab}\circ (\phi_j^{\operatorname{geo}})^{-1}\circ G_{ij}\bigr).
\]
Set
\[
K_{ab}:=\widetilde g_{ab}\circ
(\phi_j^{\operatorname{geo}})^{-1}.
\]
The Leibniz rule for three factors gives, for every multi-index $\alpha$,
\begin{equation}
\label{eq:leibniz-metrica-transicion-cap13}
D^\alpha(g_{rs}\circ\phi_i^{-1})
=
\displaystyle\sum_{\substack{\alpha_1,\alpha_2,\alpha_3\in\mathbb N_0^n\\\alpha_1+\alpha_2+\alpha_3=\alpha}}
\frac{\alpha!}{\alpha_1!\alpha_2!\alpha_3!}
D^{\alpha_1}(\partial_rG_{ij}^{a})
D^{\alpha_2}(\partial_sG_{ij}^{b})\cdot
D^{\alpha_3}(K_{ab}\circ G_{ij}).
\end{equation}
For the last factor, apply \eqref{eq:faa-di-bruno-transiciones-cap13} with $H=K_{ab}$ and with $F_{ij}$ replaced by $G_{ij}$. Its summands therefore contain derivatives $D^\delta K_{ab}$ of order $|\delta|\leq|\alpha_3|$ and products of derivatives of $G_{ij}$ of positive orders summing to $|\alpha_3|$. The first two factors of \eqref{eq:leibniz-metrica-transicion-cap13} contain derivatives of $G_{ij}$ of orders at most $|\alpha_1|+1$ and $|\alpha_2|+1$. All these quantities have uniform bounds by \ref{item:B1} and the geodesic estimates. Thus, the derivative \[
D^\alpha(g_{rs}\circ \phi_i^{-1})
\] is uniformly bounded. Therefore, for each \(k\in\mathbb{N}_{0}\), there exists a constant \(C_k'>0\) such that
\[
|D^\alpha(g_{rs}\circ \phi_i^{-1})|\le C_k',
\qquad \forall |\alpha|\le k.
\]

Now using the identity
\[
D^\alpha g_{rs}
=
\bigl(D^\alpha(g_{rs}\circ \phi_i^{-1})\bigr)\circ \phi_i,
\]
we conclude that the partial derivatives of the components \(g_{rs}\) are uniformly bounded in every chart of \(\mathcal{T}\).

To bound the inverse-metric coefficients, apply item (b) of Lemma~\ref{lema: transformacion metrica Christoffel}. For every \(p\in U_i\cap U_j^{\operatorname{geo}}\),
\[
g^{rs}(p)
=
\frac{\partial x^{r}}{\partial y^{a}}(p)
\frac{\partial x^{s}}{\partial y^{b}}(p)\,
\widetilde g^{ab}(p).
\]
Since \[
x^r
=
\pi_r\circ F_{ij}\circ \phi_j^{\operatorname{geo}}
\] on \(U_i\cap U_j^{\operatorname{geo}}\), composing with \((\phi_j^{\operatorname{geo}})^{-1}\) gives
\[
(g^{rs}\circ \phi_i^{-1}\circ F_{ij})
=
\frac{\partial F_{ij}^{r}}{\partial y^{a}}
\frac{\partial F_{ij}^{s}}{\partial y^{b}}\,
(\widetilde g^{ab}\circ (\phi_j^{\operatorname{geo}})^{-1}).
\]
Now composing with \(G_{ij}\) yields
\[
g^{rs}\circ \phi_i^{-1}
=
\left(\frac{\partial F_{ij}^{r}}{\partial y^{a}}\circ G_{ij}\right)
\left(\frac{\partial F_{ij}^{s}}{\partial y^{b}}\circ G_{ij}\right)
\bigl(\widetilde g^{ab}\circ (\phi_j^{\operatorname{geo}})^{-1}\circ G_{ij}\bigr).
\]
For this identity, the Leibniz rule again distributes $\alpha$ among three factors. If, for example,
$Q_a^r:=\partial_{y^a}F_{ij}^r$
the formula \eqref{eq:faa-di-bruno-transiciones-cap13}, applied to $Q_a^r\circ G_{ij}$, contains derivatives of $Q_a^r$ of order at most $|\alpha|$, that is, derivatives of $F_{ij}$ of order at most $|\alpha|+1$, and derivatives of $G_{ij}$ of order at most $|\alpha|$. The third factor is treated with $H=\widetilde g^{ab}\circ(\phi_j^{\operatorname{geo}})^{-1}$. The bounds in \ref{item:B1} and the geodesic bounds control all factors uniformly. We conclude that, for each multi-index \(\alpha\), there exists \(C_\alpha''>0\), independent of $i,j$, such that
\[
|D^\alpha(g^{rs}\circ \phi_i^{-1})|\le C_\alpha''.
\]
Consequently, the partial derivatives of the components \(g^{rs}\) are uniformly bounded as well.

This concludes the proof.
\end{proof}
\section{Fractional Sobolev, Besov, and Triebel--Lizorkin spaces on Riemannian manifolds of bounded geometry}
Geodesic trivializations allow us to define fractional Sobolev spaces on Riemannian manifolds of bounded geometry.
\begin{definition}
\label{def: sobolev fraccionario geometria acotada}
\index{fractional Sobolev space@fractional Sobolev space}
Let $(M,\mathbf{g})$ be a Riemannian manifold without boundary and with bounded geometry, and let \[
\mathcal T^{\operatorname{geo}}
=
(U_i^{\operatorname{geo}},
 \phi_i^{\operatorname{geo}},
 h_i^{\operatorname{geo}})_{i\in J}
\] be a geodesic trivialization as in Theorem~\ref{thm:trivializacion-geodesica}. Let $s\in\mathbb R$ and $1<p<\infty$. For $u\in\mathcal D'(M)$ and $i\in J$, the distribution $(h_i^{\operatorname{geo}}u)\circ
(\phi_i^{\operatorname{geo}})^{-1}$ has compact support in $\phi_i^{\operatorname{geo}}(U_i^{\operatorname{geo}})$. Denote its extension by zero by \[
\widetilde{
(h_i^{\operatorname{geo}}u)\circ
(\phi_i^{\operatorname{geo}})^{-1}}
\in\mathcal D'(\mathbb R^n)
\]. Define
\[
H^{s,p}(M)
:=
\left\{u\in\mathcal D'(M)
\middle|
\left(
\displaystyle\sum_{i\in J}
\left\|
\widetilde{
(h_i^{\operatorname{geo}}u)\circ
(\phi_i^{\operatorname{geo}})^{-1}}
\right\|_{H^{s,p}(\mathbb R^n)}^p
\right)^{\frac1p}<\infty
\right\}.
\]
\[
\|u\|_{H^{s,p}(M)}
:=
\left(
\displaystyle\sum_{i\in J}
\left\|
\widetilde{
(h_i^{\operatorname{geo}}u)\circ
(\phi_i^{\operatorname{geo}})^{-1}}
\right\|_{H^{s,p}(\mathbb R^n)}^p
\right)^{\frac1p}.
\]
\end{definition}

\begin{remark}
\label{obs:geometria-acotada-notemos-aunque-esta-definida-algun-abierto}
The coordinate component of $h_i^{\operatorname{geo}}u$ is understood in the sense of Proposition~\ref{prop:representacion-local-distribucion-coordenadas}. A tilde always indicates extension by zero from the coordinate open subset $\phi_i^{\operatorname{geo}}(U_i^{\operatorname{geo}})$ to $\mathbb R^n$. This extension is well defined because the support of $h_i^{\operatorname{geo}}u$ lies in $\operatorname{supp}(h_i^{\operatorname{geo}})\Subset
U_i^{\operatorname{geo}}$. Second countability of $M$ and local finiteness of the cover imply that $J$ is countable, so the sum in the definition is a well-defined series.
\end{remark}

Later we will show that this definition is independent of the chosen geodesic trivialization. To formulate this independence, it is useful to introduce the same construction for a uniformly locally finite trivialization.

\begin{definition}
\label{def:geometria-acotada-norma-avariedad-riemanniana-trivializacion-uniformemente}
\index{Sobolev space associated with an admissible trivialization@Sobolev space associated with an admissible trivialization}
Let $(M,\mathbf{g})$ be a Riemannian manifold without boundary and with bounded geometry, and let $\mathcal T=(U_i,\phi_i,h_i)_{i\in I}$ be an admissible trivialization. Let $s\in\mathbb R$ and $1<p<\infty$. For $u\in\mathcal D'(M)$, define
\[
H_{\mathcal T}^{s,p}(M)
:=
\left\{u\in\mathcal D'(M)
\middle|
\left(
\displaystyle\sum_{i\in I}
\left\|
\widetilde{(h_i u)\circ\phi_i^{-1}}
\right\|_{H^{s,p}(\mathbb R^n)}^p
\right)^{\frac1p}<\infty
\right\},
\]
with norm
\[
\|u\|_{H_{\mathcal T}^{s,p}(M)}
:=
\left(
\displaystyle\sum_{i\in I}
\left\|
\widetilde{(h_i u)\circ\phi_i^{-1}}
\right\|_{H^{s,p}(\mathbb R^n)}^p
\right)^{\frac1p}.
\]
The tilde indicates extension by zero from $\phi_i(U_i)$ to $\mathbb R^n$; it is well defined because of the auxiliary cutoff function in \ref{item:B3}.
\end{definition}
In general, the spaces $H_{\mathcal{T}}^{s,p}(M)$ depend on the chosen trivialization, and if $M$ does not have bounded geometry, they may lack desirable functional properties. For Riemannian manifolds of bounded geometry and admissible trivializations, all these Sobolev spaces coincide and have good properties, as we will see below.

\begin{lemma}[Uniform transfer between two localizations]
\label{lem:transferencia-uniforme-localizaciones}
Let $(M,\mathbf{g})$ be a Riemannian manifold without boundary and with bounded geometry. Let \[
\mathcal T=(U_i,\phi_i,h_i)_{i\in I},
\qquad
\mathcal T^\sharp
=(U_j^\sharp,\phi_j^\sharp,h_j^\sharp)_{j\in J}
\] be two admissible trivializations of $M$. Let $(\chi_i)_{i\in I}$ and $(\chi_j^\sharp)_{j\in J}$ be families of auxiliary cutoff functions, so that
\[
\operatorname{supp}(\chi_i)\Subset U_i,
\qquad
\operatorname{supp}(\chi_j^\sharp)\Subset U_j^\sharp,
\]
and each function equals one on a neighborhood of the support of the corresponding partition-of-unity element. For $U_i\cap U_j^\sharp\neq\varnothing$, set
\[
\Omega_{ij}:=\phi_i(U_i\cap U_j^\sharp),
\qquad
\Omega_{ji}^\sharp:=\phi_j^\sharp(U_i\cap U_j^\sharp)
\]
and
\[
\Phi_{ji}:=\phi_j^\sharp\circ\phi_i^{-1}\colon
\Omega_{ij}\longrightarrow\Omega_{ji}^\sharp.
\]
The function \[
\eta_{ij}:=(h_i\chi_j^\sharp)\circ\phi_i^{-1}
\] has compact support in $\Omega_{ij}$. Denote its extension by zero by $\widetilde\eta_{ij}\in C_c^\infty(\mathbb R^n)$.

For a distribution $v\in\mathcal D'(\mathbb R^n)$, the expression $v\circ\Phi_{ji}$ means the pullback of $v\restriction_{\Omega_{ji}^\sharp}$ by $\Phi_{ji}$, defined by
\begin{equation}
\label{eq:pullback-distribucional-transferencia-cap13}
\langle v\circ\Phi_{ji},\varphi\rangle_{\Omega_{ij}}
:=
\left\langle
v,
\widetilde{
(\varphi\circ\Phi_{ji}^{-1})
\left|\det D\Phi_{ji}^{-1}\right|}
\right\rangle_{\mathbb R^n},
\qquad
\varphi\in C_c^\infty(\Omega_{ij}).
\end{equation}
The tilde on the right-hand side denotes extension by zero from $\Omega_{ji}^\sharp$ to $\mathbb R^n$. This is smooth because the function being extended has compact support in $\Omega_{ji}^\sharp$. Define
\begin{equation}
\label{eq:def-operador-transferencia-localizaciones}
P_{ij}v
:=
\widetilde{\eta_{ij}(v\circ\Phi_{ji})},
\end{equation}
where the tilde on the right-hand side indicates extension by zero from $\Omega_{ij}$ to $\mathbb R^n$.

If $A^s(\mathbb R^n)$ is one of the spaces
\[
H^{s,p}(\mathbb R^n),
\qquad
F^s_{p,q}(\mathbb R^n),
\qquad
B^s_{p,q}(\mathbb R^n)
\]
in the Banach ranges of Lemma~\ref{lema: multiplicadores y difeomorfismos sobolev}, there exists a constant $C>0$, independent of $i$ and $j$, such that
\begin{equation}
\label{eq:transferencia-uniforme-localizaciones}
\|P_{ij}v\|_{A^s(\mathbb R^n)}
\leq C\|v\|_{A^s(\mathbb R^n)}.
\end{equation}
Moreover, for every $u\in\mathcal D'(M)$,
\begin{equation}
\label{eq:identidad-transferencia-localizaciones}
P_{ij}\left(
\widetilde{
(h_j^\sharp u)\circ(\phi_j^\sharp)^{-1}}
\right)
=
\widetilde{
(h_i h_j^\sharp u)\circ\phi_i^{-1}}.
\end{equation}
\end{lemma}

\begin{proof}
The support of $h_i\chi_j^\sharp$ is contained in
\[
\operatorname{supp}(h_i)\cap\operatorname{supp}(\chi_j^\sharp)
\Subset U_i\cap U_j^\sharp.
\]
Thus, $\eta_{ij}$ vanishes on a relative neighborhood of $\partial\Omega_{ij}$, and its extension by zero is smooth. On $\Omega_{ij}$,
\[
\eta_{ij}
=
(h_i\circ\phi_i^{-1})
\left[
(\chi_j^\sharp\circ(\phi_j^\sharp)^{-1})
\circ\Phi_{ji}
\right].
\]
The Leibniz rule, the bounds in \ref{item:B2} and \ref{item:B3}, and the Faà di Bruno formula of Theorem~\ref{faa di bruno multivariable} show that, for each $k\in\mathbb N_0$, there exists $C_k>0$, independent of $i$ and $j$, such that
\[
\|\widetilde\eta_{ij}\|_{BC^k(\mathbb R^n)}\leq C_k.
\]
The derivatives of $\Phi_{ji}$ and $\Phi_{ji}^{-1}$ are uniformly controlled by \ref{item:B1} and Remark~\ref{obs:geometria-acotada-item-b1-satisface-algun-atlas-geodesico}.

Let $K$ be a common bound for the entries of $D\Phi_{ji}$ and $D\Phi_{ji}^{-1}$ on all overlaps. The Leibniz formula for the determinant gives
\[
|\det D\Phi_{ji}|\leq n!K^n.
\]
Moreover,
\[
\frac{1}{|\det D\Phi_{ji}(x)|}
=
|\det D\Phi_{ji}^{-1}(\Phi_{ji}(x))|
\leq n!K^n.
\]
Thus, the Jacobians are bounded above and uniformly bounded away from zero. To apply the localized version of the coordinate-change lemma, we also verify the margins. Choose the cutoffs according to Lemma~\ref{lem:margen-uniforme-cortes-admisibles}. The support of $h_i$ lies in the region where $\chi_i=1$, and the support of $\chi_j^\sharp$ is separated from the boundary of its chart by a fixed distance. The first-derivative bounds for the charts and their inverses convert these two margins into a ball of common radius $\rho>0$ around each point of the support of $\eta_{ij}$, contained in $\Omega_{ij}$, and an analogous margin in chart $j$. To verify the first assertion, follow a coordinate segment from that point: as long as its image remains in both charts, its speed in the second is uniformly bounded; a segment shorter than the margin divided by that bound cannot reach the boundary of the second chart. The margin in the first chart ensures that the segment is defined throughout this path.

Condition \ref{item:B1} at order zero places all coordinate domains in a fixed Euclidean ball. We can cover the support of $\eta_{ij}$ by a uniformly bounded number of balls of radius $\rho/4$, with subordinate cutoffs whose derivatives have uniform bounds. On each larger ball, of radius $\rho$, the hypotheses of the localized version of Lemma~\ref{lema: multiplicadores y difeomorfismos sobolev} hold. For order $s$, fix an integer $k>|s|+n+1$; derivatives of the multipliers through order $k$ and of the transitions through order $k+1$ are controlled by \ref{item:B1}--\ref{item:B3}. The local estimate, summed over this uniformly finite cover, proves \eqref{eq:transferencia-uniforme-localizaciones} with a common constant.

Finally, $\chi_j^\sharp h_j^\sharp=h_j^\sharp$. On $\Omega_{ij}$, the definition of the pullback of distributions and compatibility of coordinate representations give
\[
\eta_{ij}
\left[
\left(
\widetilde{
(h_j^\sharp u)\circ(\phi_j^\sharp)^{-1}}
\right)\circ\Phi_{ji}
\right]
=
(h_i\chi_j^\sharp h_j^\sharp u)\circ\phi_i^{-1}.
\]
Both sides have compact support in $\Omega_{ij}$. Extending them by zero to $\mathbb R^n$ proves \eqref{eq:identidad-transferencia-localizaciones}.
\end{proof}

\begin{theorem}
\label{teo:geometria-acotada-variedad-riemanniana-geometria-acotada-trivializacion-admisible}
\index{Sobolev space@Sobolev space!independence of the admissible trivialization}
Let $(M,\mathbf{g})$ be a Riemannian manifold without boundary and with bounded geometry, and let $\mathcal T=(U_i,\phi_i,h_i)_{i\in I}$ be an admissible trivialization of $M$. Let $s\in\mathbb R$ and $1<p<\infty$. Then
\[
H_{\mathcal T}^{s,p}(M)=H^{s,p}(M),
\]
where $H^{s,p}(M)$ is defined as in Definition~\ref{def: sobolev fraccionario geometria acotada}. More precisely, there exist constants $c,C>0$, independent of $u$, such that
\[
c\|u\|_{H^{s,p}(M)}
\leq
\|u\|_{H_{\mathcal T}^{s,p}(M)}
\leq
C\|u\|_{H^{s,p}(M)}.
\]
\end{theorem}

\begin{proof}
Let \[
\mathcal T^{\operatorname{geo}}
=
(U_j^{\operatorname{geo}},\phi_j^{\operatorname{geo}},
h_j^{\operatorname{geo}})_{j\in J}
\] be the trivialization in Theorem~\ref{thm:trivializacion-geodesica}. For $i\in I$ and $j\in J$, set
\[
J_i:=\{j\in J\mid U_i\cap U_j^{\operatorname{geo}}\neq\varnothing\},
\qquad
I_j:=\{i\in I\mid U_i\cap U_j^{\operatorname{geo}}\neq\varnothing\}.
\]
The two cross-overlap bounds in \ref{item:B1}, transferred to the geodesic atlas by Remark~\ref{obs:geometria-acotada-item-b1-satisface-algun-atlas-geodesico}, provide integers $N_1,N_2$ such that
\[
\#J_i\leq N_1,
\qquad
\#I_j\leq N_2
\]
for all indices.

The locally finite identity $\displaystyle\sum_{j\in J}h_j^{\operatorname{geo}}=1$ can be multiplied by the distribution $h_i u$. After passing to chart $i$ and extending each component by zero, we obtain in $\mathcal D'(\mathbb R^n)$
\[
\widetilde{(h_i u)\circ\phi_i^{-1}}
=
\displaystyle\sum_{j\in J_i}
\widetilde{
(h_i h_j^{\operatorname{geo}}u)\circ\phi_i^{-1}}.
\]
Lemma~\ref{lem:transferencia-uniforme-localizaciones}, applied with $\mathcal T^\sharp=\mathcal T^{\operatorname{geo}}$, gives a constant $C_0>0$, independent of $i$ and $j$, such that
\[
\left\|
\widetilde{
(h_i h_j^{\operatorname{geo}}u)\circ\phi_i^{-1}}
\right\|_{H^{s,p}(\mathbb R^n)}
\leq
C_0
\left\|
\widetilde{
(h_j^{\operatorname{geo}}u)\circ
(\phi_j^{\operatorname{geo}})^{-1}}
\right\|_{H^{s,p}(\mathbb R^n)}.
\]
Set
\[
a_j
:=
\left\|
\widetilde{
(h_j^{\operatorname{geo}}u)\circ
(\phi_j^{\operatorname{geo}})^{-1}}
\right\|_{H^{s,p}(\mathbb R^n)}.
\]
Since each sum has at most $N_1$ terms,
\[
\left\|
\widetilde{(h_i u)\circ\phi_i^{-1}}
\right\|_{H^{s,p}(\mathbb R^n)}^p
\leq
C_0^pN_1^{p-1}
\displaystyle\sum_{j\in J_i}a_j^p.
\]
Summing over $i$ and rearranging the nonnegative series by Tonelli's Theorem~\ref{teo:tonelli} yields
\begin{equation*}
\begin{split}
\|u\|_{H_{\mathcal T}^{s,p}(M)}^p
&\leq
C_0^pN_1^{p-1}
\displaystyle\sum_{i\in I}\displaystyle\sum_{j\in J_i}a_j^p\\
&=
C_0^pN_1^{p-1}
\displaystyle\sum_{j\in J}\displaystyle\sum_{i\in I_j}a_j^p\\
&\leq
C_0^pN_1^{p-1}N_2
\|u\|_{H^{s,p}(M)}^p.
\end{split}
\end{equation*}
The factor $N_1$ controls the number of terms in each localization, and $N_2$ controls how many times a geodesic localization appears in the global sum.

For the reverse inequality, use $\displaystyle\sum_{i\in I}h_i=1$. In the geodesic chart $j$,
\[
\widetilde{
(h_j^{\operatorname{geo}}u)\circ
(\phi_j^{\operatorname{geo}})^{-1}}
=
\displaystyle\sum_{i\in I_j}
\widetilde{
(h_j^{\operatorname{geo}}h_i u)\circ
(\phi_j^{\operatorname{geo}})^{-1}}.
\]
Apply the same lemma with the two trivializations interchanged and repeat the preceding count. Thus,
\[
\|u\|_{H^{s,p}(M)}
\leq
C'\|u\|_{H_{\mathcal T}^{s,p}(M)}.
\]
The two estimates prove equality of the spaces as subspaces of $\mathcal D'(M)$ and, for every $u$ in that space,
\[
(C')^{-1}\|u\|_{H^{s,p}(M)}
\leq
\|u\|_{H_{\mathcal T}^{s,p}(M)}
\leq
C_0N_1^{1-\frac{1}{p}}N_2^{\frac{1}{p}}\|u\|_{H^{s,p}(M)}.
\]
\end{proof}
Having defined fractional Sobolev spaces on Riemannian manifolds, we turn to Besov and Triebel--Lizorkin spaces. On $\mathbb{R}^{n}$, priority is often given to Besov spaces, which at first sight seem simpler than Triebel--Lizorkin spaces. On a noncompact manifold, the situation changes: Besov spaces do not satisfy the strong localization principle with outer $\ell^p(I)$ summation when the fine index differs from $p$, as proved in Theorem~\ref{teo:principio-localizacion-F-H-W}; see also \cite[Theorem~2.4.7(i)]{Triebel2}. We will therefore define Triebel--Lizorkin spaces first by localization and then construct Besov spaces by real interpolation.
\begin{definition}
\label{def: Triebel-Lizorkin y Besov en variedades}
\index{Triebel--Lizorkin and Besov spaces on manifolds@Triebel--Lizorkin and Besov spaces on manifolds}
Let $(M,\mathbf{g})$ be a Riemannian manifold without boundary and with bounded geometry, let $\mathcal T=(U_i,\phi_i,h_i)_{i\in I}$ be an admissible trivialization, and let $s\in\mathbb R$.
\begin{enumerate}[label=(\alph*)]
\item If $1\leq p<\infty$ and $1\leq q\leq\infty$, define
\[
F^s_{p,q,\mathcal T}(M)
:=
\left\{u\in\mathcal D'(M)
\middle|
\|u\|_{F^s_{p,q,\mathcal T}(M)}<\infty
\right\},
\]
where
\[
\|u\|_{F^s_{p,q,\mathcal T}(M)}
:=
\left(
\displaystyle\sum_{i\in I}
\left\|
\widetilde{(h_i u)\circ\phi_i^{-1}}
\right\|_{F^s_{p,q}(\mathbb R^n)}^p
\right)^{\frac1p}.
\]
\item If $1\leq p<\infty$, $1\leq q\leq\infty$, and $-\infty<s_0<s<s_1<\infty$, write $s=(1-\Theta)s_0+\Theta s_1$, with $0<\Theta<1$, and define
\[
B^s_{p,q,\mathcal T}(M)
:=
\bigl(F^{s_0}_{p,p,\mathcal T}(M),
F^{s_1}_{p,p,\mathcal T}(M)\bigr)_{\Theta,q}.
\]
\end{enumerate}
In the first item, the tilde indicates extension by zero of the coordinate component from $\phi_i(U_i)$ to $\mathbb R^n$. In the indicated ranges, these expressions are norms.
\end{definition}

\begin{remark}\label{obs:geometria-acotada-trivializacion-geodesica-teorema-thm-trivializacion-geodesica}
\begin{enumerate}
\item If $\mathcal{T}$ is the geodesic trivialization in Theorem~\ref{thm:trivializacion-geodesica}, we omit the reference to $\mathcal{T}$ and simply write
\[
 F_{p,q}^{s}(M)
 \qquad\text{and}\qquad
 B_{p,q}^{s}(M).
 \]

\item In the range adopted here, the preceding expressions define norms. Triebel's theory also allows $0<p<1$ or $0<q<1$, but requires atomic methods for quasi-Banach spaces that will not be needed in this book.

\item The definition of $F_{p,q,\mathcal{T}}^{s}(M)$ uses coordinate localization. For Besov spaces, however, we adopt the definition by real interpolation, since the same localization principle as for Triebel--Lizorkin spaces is not generally available.

\item For every $s\in\mathbb R$ and every $1\leq p<\infty$,
\[
 B_{p,p,\mathcal{T}}^{s}(M)=F_{p,p,\mathcal{T}}^{s}(M)
 \]
and there exist constants $c_{s,p,\mathcal T},C_{s,p,\mathcal T}>0$ such that
\[
 c_{s,p,\mathcal T}\|u\|_{F^s_{p,p,\mathcal T}(M)}
 \leq\|u\|_{B^s_{p,p,\mathcal T}(M)}
 \leq C_{s,p,\mathcal T}\|u\|_{F^s_{p,p,\mathcal T}(M)}
 \]
for every $u$ in the common space.

\item The definition of $B_{p,q,\mathcal{T}}^{s}(M)$ is independent of the choice of $s_{0}$ and $s_{1}$, provided that
\[
 -\infty<s_{0}<s<s_{1}<\infty
 \qquad\text{and}\qquad
 s=(1-\Theta)s_{0}+\Theta s_{1}.
 \]
\item Alternatively, we could have defined
\[B_{p,p}^{s}(M):=\left(H^{s_{0},p}(M),H^{s_{1},p}(M)\right)_{\Theta,p},\]
where $s_{0},s_{1}\in \mathbb{R}$, $1<p<\infty$, $0<\Theta<1$, and $s=(1-\Theta)s_{0}+\Theta s_{1}$.
\end{enumerate}
\end{remark}

\begin{proof}[Justification of items \textup{(d)}--\textup{(f)}]
To apply interpolation, we construct order-independent analysis and synthesis operators defining a retraction for all the spaces under consideration. Fix the auxiliary cutoff functions $\chi_i$ from \ref{item:B3}, set $\Omega_i:=\phi_i(U_i)$, and, for $r\in\mathbb R$, write
\[
 \mathbb X_r
 :=
 \ell^p\bigl(I;F^r_{p,p}(\mathbb R^n)\bigr).
\]
The analysis operator is
\[
 \mathcal S u
 :=
 \left(
 \widetilde{(h_i u)\circ\phi_i^{-1}}
 \right)_{i\in I}.
\]
By definition, $\mathcal S\colon F^r_{p,p,\mathcal T}(M)\to\mathbb X_r$ is an isometry.

We define synthesis precisely. In chart $i$, let
\[
 a_i(x)
 :=
 \sqrt{\det\bigl(\mathbf{g}^{(i)}(\phi_i^{-1}(x))\bigr)}.
\]
If $z\in\mathcal D'(\mathbb R^n)$ has compact support contained in $\Omega_i$, the distribution $\mathcal E_i z\in\mathcal D'(M)$ is determined by
\[
 \langle\mathcal E_i z,\psi\rangle
 :=
 \left\langle
 z,
 a_i\bigl((\psi\restriction_{U_i})\circ\phi_i^{-1}\bigr)
 \right\rangle,
 \qquad \psi\in C_c^\infty(M).
\]
This expression is well defined because it uses the test function only on a neighborhood of the compact support of $z$; this is the extension of the action established in Proposition~\ref{prop: extension soporte compacto}. For $v=(v_i)_{i\in I}$, set
\[
 \mathcal Rv
 :=
 \displaystyle\sum_{i\in I}
 \mathcal E_i\bigl(\widetilde\chi_i v_i\bigr),
 \qquad
 \widetilde\chi_i
 :=
 \widetilde{\chi_i\circ\phi_i^{-1}}.
\]
The sum defines a distribution: only finitely many terms contribute on the support of a test function. Since $\chi_i=1$ near $\operatorname{supp}(h_i)$ and $\displaystyle \sum_{i\in I} h_i=1$, we obtain, first on each test function and therefore in $\mathcal D'(M)$,
\[
 \mathcal R\mathcal S u
 =
 \displaystyle\sum_{i\in I}h_i u
 =u.
\]

We show that synthesis is continuous at every order. Localizing $\mathcal Rv$ with $h_j$ involves only the indices
\[
 I_j:=\{i\in I\mid U_i\cap U_j\neq\varnothing\}.
\]
The family $I_j$ has uniformly bounded cardinality. Each summand, written in chart $j$, is obtained from $\widetilde\chi_i v_i$ by a coordinate change and multiplication by a smooth function with uniformly bounded derivatives. Lemma~\ref{lem:transferencia-uniforme-localizaciones} gives, with a constant $C_r$ independent of $i$ and $j$,
\[
 \left\|
 \widetilde{\bigl(h_j\mathcal E_i(\widetilde\chi_i v_i)\bigr)
 \circ\phi_j^{-1}}
 \right\|_{F^r_{p,p}(\mathbb R^n)}
 \leq C_r\|v_i\|_{F^r_{p,p}(\mathbb R^n)}.
\]
H\"older's inequality for the finite sum, followed by Tonelli, then shows that there exists $C'_r>0$ such that
\[
 \|\mathcal Rv\|_{F^r_{p,p,\mathcal T}(M)}
 \leq C'_r\|v\|_{\mathbb X_r}.
\]
Thus, $F^r_{p,p,\mathcal T}(M)$ is a retract of $\mathbb X_r$, and the operators $\mathcal S$ and $\mathcal R$ do not depend on $r$. Moreover, $\mathcal S\mathcal R$ is a continuous projection whose closed image is $\mathcal S(F^r_{p,p,\mathcal T}(M))$. Since $\mathbb X_r$ is Banach and $\mathcal S$ is isometric, $F^r_{p,p,\mathcal T}(M)$ is Banach as well. This verifies that the pairs interpolated below are Banach pairs.

Now let $r_0<r_1$, $0<\eta<1$, and $r=(1-\eta)r_0+\eta r_1$. The Euclidean identity
\[
 (F^{r_0}_{p,p}(\mathbb R^n),F^{r_1}_{p,p}(\mathbb R^n))_{\eta,p}
 =F^r_{p,p}(\mathbb R^n)
\]
follows from the dyadic decomposition, $B^t_{p,p}(\mathbb R^n)=F^t_{p,p}(\mathbb R^n)$, and Proposition~\ref{prop:interpolacion-real-ellp}; it also holds for $p=1$. The same proposition gives
\[
 (\mathbb X_{r_0},\mathbb X_{r_1})_{\eta,p}=\mathbb X_r.
\]
Theorem~\ref{teo:interpolacion-operadores-metodo-K}, applied to $\mathcal S$ and $\mathcal R$, preserves the identity $\mathcal R\mathcal S=I$ on the interpolated spaces. Consequently,
\begin{equation}
 \label{eq:interpolacion-Fpp-localizado-cap13}
 (F^{r_0}_{p,p,\mathcal T}(M),F^{r_1}_{p,p,\mathcal T}(M))_{\eta,p}
 =F^r_{p,p,\mathcal T}(M)
\end{equation}
with equivalence of norms. Taking $r_0=s_0$, $r_1=s_1$, and $\eta=\Theta$ proves item \textup{(d)}.

For item \textup{(e)}, consider two admissible pairs of endpoints. Choose $a<s<b$ so that all four endpoints lie in $(a,b)$. By \eqref{eq:interpolacion-Fpp-localizado-cap13}, each endpoint space is a space $(F^a_{p,p,\mathcal T}(M),F^b_{p,p,\mathcal T}(M))_{\rho,p}$, with parameter $\rho$ determined affinely by its order. The real reiteration Theorem~\ref{teo:reiteracion-interpolacion-real} reduces both definitions of $B^s_{p,q,\mathcal T}(M)$ to the same space
\[
 (F^a_{p,p,\mathcal T}(M),F^b_{p,p,\mathcal T}(M))_{\rho_s,q},
 \qquad
 \rho_s=\frac{s-a}{b-a}.
\]
This proves independence of the endpoints, including equivalence of the resulting norms.

Finally, if $1<p<\infty$, the same analysis and synthesis are continuous when $F^r_{p,p}(\mathbb R^n)$ is replaced by $H^{r,p}(\mathbb R^n)$; the proof uses exactly the preceding multiplication, coordinate-change, and overlap estimates. Item \textup{(a)} of Theorem~\ref{teo:interpolacion-escalas-BF} identifies real interpolation of the Euclidean Bessel spaces with $B^s_{p,p}(\mathbb R^n)=F^s_{p,p}(\mathbb R^n)$. Interpolating the retraction gives
\[
 (H^{s_0,p}(M),H^{s_1,p}(M))_{\Theta,p}
 =F^s_{p,p}(M)=B^s_{p,p}(M),
\]
which is item \textup{(f)}.
\end{proof}

These spaces are independent of the chosen trivialization (since the Riemannian manifold has bounded geometry), as we now show:
\begin{theorem}
\label{teo:geometria-acotada-variedad-riemanniana-trivializacion-admisible}
\index{Triebel--Lizorkin and Besov spaces@Triebel--Lizorkin and Besov spaces!independence of the admissible trivialization}
Let $(M,\mathbf{g})$ be a Riemannian manifold without boundary and with bounded geometry, and let $\mathcal T=(U_i,\phi_i,h_i)_{i\in I}$ be an admissible trivialization. Let $s\in\mathbb R$, $1\leq p<\infty$, and $1\leq q\leq\infty$. Then
\[
F^s_{p,q,\mathcal T}(M)=F^s_{p,q}(M)
\qquad\text{and}\qquad
B^s_{p,q,\mathcal T}(M)=B^s_{p,q}(M)
\]
and there exist positive constants $c_F,C_F,c_B,C_B$, depending on $s,p,q$ and the uniform data of the two trivializations but not on $u$, such that
\[
 c_F\|u\|_{F^s_{p,q}(M)}
 \leq\|u\|_{F^s_{p,q,\mathcal T}(M)}
 \leq C_F\|u\|_{F^s_{p,q}(M)}
\]
and
\[
 c_B\|u\|_{B^s_{p,q}(M)}
 \leq\|u\|_{B^s_{p,q,\mathcal T}(M)}
 \leq C_B\|u\|_{B^s_{p,q}(M)}.
\]
\end{theorem}

\begin{proof}
For the scale $F$, repeat the count in the proof of Theorem~\ref{teo:geometria-acotada-variedad-riemanniana-geometria-acotada-trivializacion-admisible}. If $\mathcal T^{\operatorname{geo}}$ is the geodesic trivialization, in $\mathcal D'(\mathbb R^n)$ we have
\[
\widetilde{(h_i u)\circ\phi_i^{-1}}
=
\displaystyle\sum_{j\in J_i}
\widetilde{
(h_i h_j^{\operatorname{geo}}u)\circ\phi_i^{-1}}.
\]
Lemma~\ref{lem:transferencia-uniforme-localizaciones}, with $A^s=F^s_{p,q}$, gives
\[
\left\|
\widetilde{
(h_i h_j^{\operatorname{geo}}u)\circ\phi_i^{-1}}
\right\|_{F^s_{p,q}(\mathbb R^n)}
\leq
C_0
\left\|
\widetilde{
(h_j^{\operatorname{geo}}u)\circ
(\phi_j^{\operatorname{geo}})^{-1}}
\right\|_{F^s_{p,q}(\mathbb R^n)}.
\]
The constant $C_0$ depends on $s,p,q$ and the uniform bounds for the cutoff functions, coordinate transitions, and their Jacobians, but not on $i$ or $j$. If $N_1,N_2$ are the two cross multiplicities, outer $\ell^p(I)$ summation gives
\[
\|u\|_{F^s_{p,q,\mathcal T}(M)}^p
\leq
C_0^pN_1^{p-1}N_2
\|u\|_{F^s_{p,q}(M)}^p.
\]
Interchanging the trivializations gives the reverse estimate. The argument remains valid when $q=\infty$, since it applies to the local norms, which are nonnegative numbers, while the outer sum retains exponent $p<\infty$.

For the scale $B$, fix $s_0<s<s_1$ and $s=(1-\theta)s_0+\theta s_1$. The identities on the common space $\mathcal D'(M)$ define compatible isomorphisms
\[
F^{s_\nu}_{p,p,\mathcal T}(M)
\longrightarrow
F^{s_\nu}_{p,p}(M),
\qquad \nu\in\{0,1\}.
\]
The operator interpolation theorem for the $K$ method~\ref{teo:interpolacion-operadores-metodo-K}, applied also to the inverse identity, gives
\[
\bigl(F^{s_0}_{p,p,\mathcal T}(M),
F^{s_1}_{p,p,\mathcal T}(M)\bigr)_{\theta,q}
=
\bigl(F^{s_0}_{p,p}(M),F^{s_1}_{p,p}(M)\bigr)_{\theta,q}
\]
with norms controlled in both directions by the norms of the identities at the two endpoints and their inverses. These are precisely the constants $c_B,C_B$ in the statement. This establishes the asserted equality of Besov spaces. The reiteration Theorem~\ref{teo:reiteracion-interpolacion-real} also shows that the construction is independent of the auxiliary endpoints.
\end{proof}
\subsection{Comparison with intrinsic definitions}

\begin{remark}[Heat on several components]
\label{obs:calor-componentes-geometria-acotada}
In this chapter, $M$ is not assumed connected. Its components $(M_\nu)_{\nu\in A}$ form an at most countable family, since they are open and $M$ has a countable base. Each is complete by Remark~\ref{obs: geometria acotada implica completa}. To apply the constructions of the heat chapter, use
\[
 L^p(M,\mathbf E)
 =\ell^p\bigl(A;L^p(M_\nu,\mathbf E|_{M_\nu})\bigr),
 \qquad 1\leq p<\infty,
\]
and define the semigroup componentwise. It is contractive because every component semigroup is contractive. If $\mathbf u=(\mathbf u_\nu)_{\nu\in A}$, then
\[
 \|P_t^{\mathbf E}\mathbf u-\mathbf u\|_{L^{p}(M,\mathbf E)}^p
 =\displaystyle\sum_{\nu\in A}
   \|P_t^{\mathbf E|_{M_\nu}}\mathbf u_\nu-\mathbf u_\nu\|_{L^{p}(M_\nu,\mathbf E\restriction_{M_\nu})}^p
 \longrightarrow0:
\]
each summand tends to zero and is bounded by $2^p\|\mathbf u_\nu\|_{L^{p}(M_\nu,\mathbf E\restriction_{M_\nu})}^p$, a convergent series. Thus, strong continuity is preserved, including for $p=1$. The projections onto the components are continuous and commute with the Bochner integrals defining the potentials. The intrinsic scales are therefore constructed with the same semigroup componentwise. The following estimates use constants independent of $\nu$, since the geometric bounds are uniform on all of $M$. In the scalar case, use the trivial bundle.
\end{remark}

The comparison will rely on three distinct mechanisms. The first is uniform transfer between charts, already established in Lemma~\ref{lem:transferencia-uniforme-localizaciones}. The second compares partial and covariant derivatives at integer orders. The third interpolates these integer-order identifications and compares them with the scale generated by the Laplacian. We begin by specifying the role of quadratic partitions.

\begin{lemma}[Equivalence of quadratic localization]
\label{lem:equivalencia-localizacion-cuadratica}
Let $(M,\mathbf{g})$ be a Riemannian manifold without boundary and with bounded geometry, let $\mathcal Q=(U_i,\phi_i,h_i)_{i\in I}$ be an admissible quadratic localization system, let $\rho_i:=h_i^2$, and let $\mathcal T_{\mathcal Q}=(U_i,\phi_i,\rho_i)_{i\in I}$ be the associated trivialization. For $s\in\mathbb R$, $1<p<\infty$, and $1\leq q\leq\infty$, there exist constants $c_H,C_H,c_F,C_F>0$, independent of $u$, such that
\[
c_H\|u\|_{H^{s,p}_{\mathcal T_{\mathcal Q}}(M)}
\leq
\left(
\displaystyle\sum_{i\in I}
\left\|
\widetilde{(h_i u)\circ\phi_i^{-1}}
\right\|_{H^{s,p}(\mathbb R^n)}^p
\right)^{\frac1p}
\leq
C_H\|u\|_{H^{s,p}_{\mathcal T_{\mathcal Q}}(M)}
\]
and
\[
c_F\|u\|_{F^s_{p,q,\mathcal T_{\mathcal Q}}(M)}
\leq
\left(
\displaystyle\sum_{i\in I}
\left\|
\widetilde{(h_i u)\circ\phi_i^{-1}}
\right\|_{F^s_{p,q}(\mathbb R^n)}^p
\right)^{\frac1p}
\leq
C_F\|u\|_{F^s_{p,q,\mathcal T_{\mathcal Q}}(M)}.
\]
When $q=p$, the second comparison also provides localization of $B^s_{p,p}(M)$.
\end{lemma}

\begin{proof}
Fix one of the two choices
\[
A^s(\mathbb R^n)=H^{s,p}(\mathbb R^n)
\qquad\text{or}\qquad
A^s(\mathbb R^n)=F^s_{p,q}(\mathbb R^n)
\]
and set
\[
N_h(u)
:=
\left(
\displaystyle\sum_{i\in I}
\left\|
\widetilde{(h_i u)\circ\phi_i^{-1}}
\right\|_{A^s(\mathbb R^n)}^p
\right)^{\frac1p},
\]
\[
N_\rho(u)
:=
\left(
\displaystyle\sum_{i\in I}
\left\|
\widetilde{(\rho_i u)\circ\phi_i^{-1}}
\right\|_{A^s(\mathbb R^n)}^p
\right)^{\frac1p}.
\]
By definition, $N_\rho$ is the norm associated with $\mathcal T_{\mathcal Q}$. Since $\operatorname{supp}(h_i)\Subset U_i$, the function $h_i\circ\phi_i^{-1}$ has an extension by zero $\widetilde{h_i\circ\phi_i^{-1}}\in C_c^\infty(\mathbb R^n)$, and
\[
\widetilde{(\rho_i u)\circ\phi_i^{-1}}
=
\widetilde{h_i\circ\phi_i^{-1}}\,
\widetilde{(h_i u)\circ\phi_i^{-1}}.
\]
The multiplier theorem in Lemma~\ref{lema: multiplicadores y difeomorfismos sobolev} and the uniform bounds on \ref{item:B2} give a constant $K_0>0$, independent of $i\in I$, such that
\[
\left\|
\widetilde{(\rho_i u)\circ\phi_i^{-1}}
\right\|_{A^s(\mathbb R^n)}
\leq
K_0
\left\|
\widetilde{(h_i u)\circ\phi_i^{-1}}
\right\|_{A^s(\mathbb R^n)}.
\]
Raising to the power $p$ and summing over $i\in I$ gives
\begin{equation}
\label{eq:cuadratica-primera-direccion}
N_\rho(u)\leq K_0N_h(u).
\end{equation}

For the reverse estimate, use the locally finite identity
\[
\displaystyle\sum_{j\in I}\rho_j=1.
\]
For each $i\in I$, define
\[
I_i:=\{j\in I\mid U_i\cap U_j\neq\varnothing\}.
\]
Uniform finiteness of the cover provides $L\in\mathbb N$ such that $\#I_i\leq L$ for every $i\in I$. Therefore,
\begin{equation}
\label{eq:descomposicion-cuadratica-local-cap13}
\widetilde{(h_i u)\circ\phi_i^{-1}}
=
\displaystyle\sum_{j\in I_i}
\widetilde{(h_i\rho_j u)\circ\phi_i^{-1}}
\quad\text{on }\mathcal D'(\mathbb R^n).
\end{equation}

Fix $j\in I_i$ and let $\chi_j$ be the auxiliary cutoff function associated with $\rho_j$. On the overlap, set
\[
\Omega_{ij}:=\phi_i(U_i\cap U_j),
\qquad
\Phi_{ji}:=\phi_j\circ\phi_i^{-1}.
\]
In accordance with \eqref{eq:pullback-distribucional-transferencia-cap13}, define
\[
P_{ij}^{h,\rho}z
:=
\widetilde{
\bigl((h_i\chi_j)\circ\phi_i^{-1}\bigr)
(z\circ\Phi_{ji})}.
\]
The support of $(h_i\chi_j)\circ\phi_i^{-1}$ is compact in $\Omega_{ij}$. The proof of Lemma~\ref{lem:transferencia-uniforme-localizaciones}, which does not use the partition-of-unity identity for this estimate, provides a constant $K_1>0$, independent of $i$ and $j$, such that
\begin{equation}
\label{eq:cota-transferencia-cuadratica-cap13}
\|P_{ij}^{h,\rho}z\|_{A^s(\mathbb R^n)}
\leq
K_1\|z\|_{A^s(\mathbb R^n)}.
\end{equation}
Since $\chi_j\rho_j=\rho_j$,
\[
P_{ij}^{h,\rho}
\left(
\widetilde{(\rho_j u)\circ\phi_j^{-1}}
\right)
=
\widetilde{(h_i\rho_j u)\circ\phi_i^{-1}}.
\]
The preceding identities imply
\[
\left\|
\widetilde{(h_i u)\circ\phi_i^{-1}}
\right\|_{A^s(\mathbb R^n)}
\leq
K_1\displaystyle\sum_{j\in I_i}
\left\|
\widetilde{(\rho_j u)\circ\phi_j^{-1}}
\right\|_{A^s(\mathbb R^n)}.
\]
Hölder's inequality for a sum of at most $L$ terms gives
\[
\left(\displaystyle\sum_{j\in I_i}a_j\right)^p
\leq
L^{p-1}\displaystyle\sum_{j\in I_i}a_j^p
\]
for every family of nonnegative numbers $(a_j)_{j\in I_i}$. Applying this inequality and summing over $i\in I$ gives
\begin{equation*}
\begin{split}
N_h(u)^p
&\leq
K_1^pL^{p-1}
\displaystyle\sum_{i\in I}\displaystyle\sum_{j\in I_i}
\left\|
\widetilde{(\rho_j u)\circ\phi_j^{-1}}
\right\|_{A^s(\mathbb R^n)}^p\\
&=
K_1^pL^{p-1}
\displaystyle\sum_{j\in I}\#I_j
\left\|
\widetilde{(\rho_j u)\circ\phi_j^{-1}}
\right\|_{A^s(\mathbb R^n)}^p\\
&\leq
K_1^pL^pN_\rho(u)^p.
\end{split}
\end{equation*}
The intermediate equality uses symmetry of the relation $U_i\cap U_j\neq\varnothing$. Consequently,
\begin{equation}
\label{eq:cuadratica-segunda-direccion}
N_h(u)\leq K_1L\,N_\rho(u).
\end{equation}
Estimates \eqref{eq:cuadratica-primera-direccion} and \eqref{eq:cuadratica-segunda-direccion} prove the two comparisons. For $A^s=F^s_{p,p}$, also use the Euclidean identification \[
F^s_{p,p}(\mathbb R^n)=B^s_{p,p}(\mathbb R^n)
\] of Corollary~\ref{teo:identificaciones-H-W-B-F}.
\end{proof}

\begin{proposition}[Retraction, interpolation, and distributional duality]
\label{prop:retraccion-interpolacion-dualidad-geometria-acotada}
Let $(M,\mathbf{g})$ be a Riemannian manifold without boundary and with bounded geometry, $s_0,s_1\in\mathbb R$, $s_0<s_1$, $1<p<\infty$, and $0<\theta<1$. If $s=(1-\theta)s_0+\theta s_1$, then
\begin{equation}
\label{eq:interpolacion-compleja-H-localizado}
[H^{s_0,p}(M),H^{s_1,p}(M)]_\theta
=
H^{s,p}(M)
\end{equation}
and there exist constants $c_I,C_I>0$, depending on $s_0,s_1,\theta,p$ and the uniform localization data, such that
\[
 c_I\|u\|_{H^{s,p}(M)}
 \leq
 \|u\|_{[H^{s_0,p}(M),H^{s_1,p}(M)]_\theta}
 \leq C_I\|u\|_{H^{s,p}(M)}.
\]

Moreover, the identification of locally integrable functions with regular distributions, constructed using $d\lambda_{\mathbf{g}}$, induces on $C_c^\infty(M)$ the bilinear pairing
\begin{equation}
\label{eq:emparejamiento-distribucional-riemann-lebesgue}
\langle v,u\rangle_{\mathbf{g}}
:=
T_v(u)
=
\int_Muv\,d\lambda_{\mathbf{g}},
\qquad u,v\in C_c^\infty(M).
\end{equation}
For every $r\in\mathbb R$, this pairing has a unique continuous bilinear extension
\[
\langle\cdot,\cdot\rangle_{r,p;\mathbf{g}}\colon
H^{-r,p'}(M)\times H^{r,p}(M)
\longrightarrow
\mathbb K,
\]
and the operator
\[
\mathcal J_{r,p;\mathbf{g}}\colon
H^{-r,p'}(M)
\longrightarrow
\bigl(H^{r,p}(M)\bigr)',
\qquad
\mathcal J_{r,p;\mathbf{g}}(v)(u)
:=
\langle v,u\rangle_{r,p;\mathbf{g}},
\]
is a topological isomorphism. In particular,
\begin{equation}
\label{eq:dualidad-H-localizado}
\bigl(H^{r,p}(M)\bigr)'
\cong
H^{-r,p'}(M)
\end{equation}
and there exist constants $c_D,C_D>0$, depending on $r,p$ and the uniform localization data, such that
\[
 c_D\|v\|_{H^{-r,p'}(M)}
 \leq\|\mathcal J_{r,p;\mathbf{g}}v\|_{(H^{r,p}(M))'}
 \leq C_D\|v\|_{H^{-r,p'}(M)}.
\]
When $\mathbb K=\mathbb C$, the pairing \eqref{eq:emparejamiento-distribucional-riemann-lebesgue} involves no conjugation: it is the distributional duality induced by the Riemann--Lebesgue measure, not the Hermitian inner product of $L^2(M)$.
\end{proposition}

\begin{proof}
Start with a geodesic trivialization and apply Lemma~\ref{lem:normalizacion-cuadratica-uniforme}. This gives an admissible quadratic system \[
\mathcal Q=(U_i,\phi_i,h_i)_{i\in I},
\qquad
\displaystyle\sum_{i\in I}h_i^2=1,
\] and auxiliary cutoff functions $(\chi_i)_{i\in I}$ such that
\[
\chi_i=1
\quad\text{on a neighborhood of }\operatorname{supp}(h_i),
\qquad
\operatorname{supp}(\chi_i)\Subset U_i.
\]
Write
\[
\Omega_i:=\phi_i(U_i),
\qquad
a_i(x):=
\sqrt{\det\bigl(\mathbf{g}^{(i)}(\phi_i^{-1}(x))\bigr)},
\qquad x\in\Omega_i.
\]
The function $\chi_i\circ\phi_i^{-1}$ has compact support in $\Omega_i$. Denote by \[
\widetilde\chi_i
:=
\widetilde{\chi_i\circ\phi_i^{-1}}
\in C_c^\infty(\mathbb R^n)
\] its extension by zero from $\Omega_i$ to $\mathbb R^n$.

Before defining synthesis, we specify how to reconstruct a distribution from a chart. If $z\in\mathcal D'(\mathbb R^n)$ and $\operatorname{supp}(z)\Subset\Omega_i$, define $\mathcal E_i z\in\mathcal D'(M)$ by
\begin{equation}
\label{eq:reconstruccion-distribucion-carta-cap13}
\langle\mathcal E_i z,\psi\rangle
:=
\left\langle
z,
a_i\bigl((\psi\restriction_{U_i})\circ\phi_i^{-1}\bigr)
\right\rangle_{\mathbb R^n},
\qquad
\psi\in C_c^\infty(M).
\end{equation}
Although the function in the second argument need not have compact support in $\Omega_i$, only its restriction to a neighborhood of $\operatorname{supp}(z)$ is used. The expression is well defined by the extension of the action of a compactly supported distribution established in Proposition~\ref{prop: extension soporte compacto}. The definition of coordinate components of a distribution in Proposition~\ref{prop:representacion-local-distribucion-coordenadas} implies
\begin{equation}
\label{eq:reconstruccion-invierte-coordenadas-cap13}
\bigl((\mathcal E_i z)\restriction_{U_i}\bigr)\circ\phi_i^{-1}
=
z\restriction_{\Omega_i}
\quad\text{on }\mathcal D'(\Omega_i).
\end{equation}
Indeed, for $\varphi\in C_c^\infty(\Omega_i)$,
\[
\left\langle
\bigl((\mathcal E_i z)\restriction_{U_i}\bigr)\circ\phi_i^{-1},
\varphi
\right\rangle
=
\left\langle
\mathcal E_i z,
\left(\frac{\varphi}{a_i}\circ\phi_i\right)
\right\rangle
=
\langle z,\varphi\rangle.
\]
If $z$ is a locally integrable function with compact support in $\Omega_i$, then $\mathcal E_i z$ is the regular distribution associated with the function $z\circ\phi_i$ on $U_i$, extended by zero outside $U_i$. This construction avoids introducing notation for transporting distributions before specifying its meaning.

Define the analysis operator by
\[
\mathcal Su
:=
\left(
\widetilde{(h_i u)\circ\phi_i^{-1}}
\right)_{i\in I}
\]
and, for $v=(v_i)_{i\in I}$, the synthesis operator by
\begin{equation}
\label{eq:operador-sintesis-explicito-cap13}
\mathcal Rv
:=
\displaystyle\sum_{i\in I}
 h_i\,\mathcal E_i(\widetilde\chi_i v_i).
\end{equation}
The product $\widetilde\chi_i v_i$ has support contained in $\operatorname{supp}(\widetilde\chi_i)\Subset\Omega_i$, so $\mathcal E_i(\widetilde\chi_i v_i)$ is defined by \eqref{eq:reconstruccion-distribucion-carta-cap13}. The sum in \eqref{eq:operador-sintesis-explicito-cap13} is locally finite, since its $i$th summand has support contained in $\operatorname{supp}(h_i)$ and the family $(U_i)_{i\in I}$ is locally finite.

Definition~\ref{def:sumas-ellq-familia-banach} specifies the sequence spaces we will use. Lemma~\ref{lem:equivalencia-localizacion-cuadratica}, together with independence of the trivialization, shows that, for each $r\in\mathbb R$,
\[
\mathcal S\colon
H^{r,p}(M)
\longrightarrow
\ell^p\bigl(I;H^{r,p}(\mathbb R^n)\bigr)
\]
is continuous.

We prove continuity of $\mathcal R$. For $j\in I$, set
\[
I_j
:=
\{i\in I\mid U_i\cap U_j\neq\varnothing\}.
\]
There exists $L\geq1$ such that $\#I_j\leq L$ for every $j\in I$. If $i\in I_j$, write
\[
\Omega_{ji}:=\phi_j(U_j\cap U_i),
\qquad
\Phi_{ij}:=\phi_i\circ\phi_j^{-1}\colon
\Omega_{ji}\longrightarrow\phi_i(U_i\cap U_j),
\]
and define
\[
\eta_{ji}:=(h_jh_i\chi_i)\circ\phi_j^{-1}
\in C_c^\infty(\Omega_{ji}).
\]
Its support is compact in $\Omega_{ji}$ because
\[
\operatorname{supp}(h_jh_i\chi_i)
\subseteq
\operatorname{supp}(h_j)\cap\operatorname{supp}(h_i)
\Subset U_j\cap U_i.
\]
Let $z\in\mathcal D'(\mathbb R^n)$. On the overlap, use the restriction of $z$ to $\phi_i(U_i\cap U_j)$ and define its distributional composition with $\Phi_{ij}$ by
\begin{equation}
\label{eq:composicion-distribucional-cambio-carta-cap13}
\left\langle z\circ\Phi_{ij},\varphi\right\rangle_{\Omega_{ji}}
:=
\left\langle
z,
\widetilde{
(\varphi\circ\Phi_{ij}^{-1})
\left|\det D\Phi_{ij}^{-1}\right|}
\right\rangle_{\mathbb R^n},
\qquad
\varphi\in C_c^\infty(\Omega_{ji}).
\end{equation}
Here the tilde indicates extension by zero from $\phi_i(U_i\cap U_j)$ to $\mathbb R^n$; the extension is smooth because the function's support is compactly contained in that open set. For a locally integrable function, this definition agrees with ordinary composition $x\mapsto z(\Phi_{ij}(x))$. Set
\[
P_{ji}z
:=
\widetilde{\eta_{ji}(z\circ\Phi_{ij})},
\]
where the tilde indicates extension by zero from $\Omega_{ji}$ to $\mathbb R^n$. Compactness of the support of $\eta_{ji}$ inside $\Omega_{ji}$ guarantees that this extension is a well-defined distribution. The transformation of the volume densities on the overlap is
\begin{equation}
\label{eq:transformacion-densidad-volumen-sintesis-cap13}
a_j(x)
=
a_i(\Phi_{ij}(x))
\left|\det D\Phi_{ij}(x)\right|,
\qquad x\in\Omega_{ji}.
\end{equation}
Evaluating both sides against a test function on $C_c^\infty(\Omega_{ji})$, definition \eqref{eq:reconstruccion-distribucion-carta-cap13} and \eqref{eq:transformacion-densidad-volumen-sintesis-cap13} give
\begin{equation}
\label{eq:formula-local-sintesis-explicita-cap13}
\widetilde{
\bigl(
 h_jh_i\,\mathcal E_i(\widetilde\chi_i v_i)
\bigr)\circ\phi_j^{-1}}
=
P_{ji}v_i.
\end{equation}
The bounds in \ref{item:B1}--\ref{item:B3}, the multiplier and coordinate-change theorem in Lemma~\ref{lema: multiplicadores y difeomorfismos sobolev}, and the uniform Jacobian bounds obtained in Lemma~\ref{lem:transferencia-uniforme-localizaciones} provide a constant $C_r>0$, independent of $i$ and $j$, such that
\begin{equation}
\label{eq:cota-operador-transferencia-sintesis-cap13}
\|P_{ji}z\|_{H^{r,p}(\mathbb R^n)}
\leq
C_r\|z\|_{H^{r,p}(\mathbb R^n)}.
\end{equation}
By \eqref{eq:formula-local-sintesis-explicita-cap13},
\[
\left\|
\widetilde{(h_j\mathcal Rv)\circ\phi_j^{-1}}
\right\|_{H^{r,p}(\mathbb R^n)}
\leq
C_r\displaystyle\sum_{i\in I_j}\|v_i\|_{H^{r,p}(\mathbb R^n)}.
\]
Raise to the power $p$ and use the fact that each sum contains at most $L$ terms. Summing over $j\in I$ and counting overlaps again gives
\begin{equation*}
\begin{split}
&\displaystyle\sum_{j\in I}
\left\|
\widetilde{(h_j\mathcal Rv)\circ\phi_j^{-1}}
\right\|_{H^{r,p}(\mathbb R^n)}^p\\
&\quad\leq
C_r^pL^{p-1}
\displaystyle\sum_{j\in I}\displaystyle\sum_{i\in I_j}
\|v_i\|_{H^{r,p}(\mathbb R^n)}^p\\
&\quad=
C_r^pL^{p-1}
\displaystyle\sum_{i\in I}\#I_i\,
\|v_i\|_{H^{r,p}(\mathbb R^n)}^p\\
&\quad\leq
C_r^pL^p
\displaystyle\sum_{i\in I}\|v_i\|_{H^{r,p}(\mathbb R^n)}^p.
\end{split}
\end{equation*}
The two inequalities of Lemma~\ref{lem:equivalencia-localizacion-cuadratica} compare the left-hand side with the norm of $H^{r,p}(M)$ through $c_H$ and $C_H$. Consequently,
\begin{equation}
\label{eq:sintesis-continua-cap13}
\|\mathcal Rv\|_{H^{r,p}(M)}
\leq
C_r'\|v\|_{\ell^p(I;H^{r,p}(\mathbb R^n))}.
\end{equation}

Since $\chi_i=1$ on a neighborhood of $\operatorname{supp}(h_i)$,
\[
\widetilde\chi_i
\widetilde{(h_i u)\circ\phi_i^{-1}}
=
\widetilde{(h_i u)\circ\phi_i^{-1}}.
\]
The identity \eqref{eq:reconstruccion-invierte-coordenadas-cap13}, applied to the distribution $h_i u$, and the equality $\displaystyle\sum_{i\in I}h_i^2=1$ give
\[
\mathcal R\mathcal Su
=
\displaystyle\sum_{i\in I}h_i^2u
=u.
\]
Thus, $\mathcal S$ is a coretraction and $\mathcal R$ a retraction, and the same pair works at every real order. These identities also establish completeness before interpolation: the operator $\mathcal S\mathcal R$ is a continuous projection on the Banach space $\ell^p(I;H^{r,p}(\mathbb R^n))$, and its image is exactly $\mathcal S(H^{r,p}(M))$. This image is closed, being the kernel of $I-\mathcal S\mathcal R$. The norm comparisons for $\mathcal S$ and its inverse $\mathcal R$ identify $H^{r,p}(M)$ with that image; therefore, $H^{r,p}(M)$ is a Banach space.

This pair also proves density of $C_c^\infty(M)$ in $H^{r,p}(M)$. Finitely supported sequences with components in $C_c^\infty(\mathbb R^n)$ are dense in $\ell^p(I;H^{r,p}(\mathbb R^n))$. If a sequence of this kind approximates $\mathcal Su$, its image under $\mathcal R$ belongs to $C_c^\infty(M)$ and converges to $u$ in $H^{r,p}(M)$ by \eqref{eq:sintesis-continua-cap13}.

Apply complex interpolation. Proposition~\ref{prop:interpolacion-ellp-valores-banach} and the Euclidean identity of Theorem~\ref{teo:interpolacion-escalas-BF} give
\[
\begin{split}
&[\ell^p(I;H^{s_0,p}(\mathbb R^n)),
  \ell^p(I;H^{s_1,p}(\mathbb R^n))]_\theta\\
&\hspace{30mm}
=
\ell^p(I;H^{s,p}(\mathbb R^n)).
\end{split}
\]
Interpolate $\mathcal S$ and $\mathcal R$ at the endpoints and use $\mathcal R\mathcal S=I$. The endpoint norms of these operators, their interpolated norms, and the localization bounds produce the constants $c_I,C_I>0$ in the statement; in particular, for every $u$ in the common space of \eqref{eq:interpolacion-compleja-H-localizado},
\[
c_I\|u\|_{H^{s,p}(M)}
\leq
\|u\|_{[H^{s_0,p}(M),H^{s_1,p}(M)]_\theta}
\leq
C_I\|u\|_{H^{s,p}(M)}.
\]

We turn to duality. According to Definition~\ref{def: distribucion regular}, each $v\in L^1_{\operatorname{loc}}(M)$ determines the regular distribution
\[
T_v(u)=\int_Muv\,d\lambda_{\mathbf{g}},
\qquad u\in C_c^\infty(M).
\]
Thus, \eqref{eq:emparejamiento-distribucional-riemann-lebesgue} is the restriction to smooth functions of the duality between $\mathcal D'(M)$ and $\mathcal D(M)$ induced by the Riemann--Lebesgue measure.

The metric criterion of Lemma~\ref{lem:geometria-acotada-variedad-riemanniana-trivializacion-uniformemente-localmente-finit} and the Faà di Bruno formula of Theorem~\ref{faa di bruno multivariable} show that $a_i$, $a_i^{-1}$, and all their derivatives are bounded by constants independent of $i$. If $v\in\mathcal D'(M)$ and $u\in C_c^\infty(M)$, then
\begin{equation}
\label{eq:accion-distribucional-coordenadas-cap13}
\begin{split}
v(u)
&=
v\left(\displaystyle\sum_{i\in I}h_i^2u\right)\\
&=
\displaystyle\sum_{i\in I}
\left\langle
\widetilde{(h_i v)\circ\phi_i^{-1}},
\widetilde{
a_i\bigl((h_i u)\circ\phi_i^{-1}\bigr)}
\right\rangle_{\mathbb R^n}.
\end{split}
\end{equation}
The sum is finite because $u$ has compact support. The second equality is precisely the definition of the coordinate component of a distribution in Proposition~\ref{prop:representacion-local-distribucion-coordenadas}. When $v$ is smooth, apply the change-of-variables formula in each chart:
\[
\int_{U_i}\eta\,d\lambda_{\mathbf{g}}
=
\int_{\Omega_i}(\eta\circ\phi_i^{-1})(x)a_i(x)\,d\lambda_n(x).
\]
The preceding identity then becomes
\begin{equation}
\label{eq:emparejamiento-coordenadas-cap13}
\int_Muv\,d\lambda_{\mathbf{g}}
=
\displaystyle\sum_{i\in I}
\int_{\Omega_i}
\bigl((h_i u)\circ\phi_i^{-1}\bigr)
\bigl((h_i v)\circ\phi_i^{-1}\bigr)
a_i\,dx.
\end{equation}

Let $v\in H^{-r,p'}(M)$ and $u\in C_c^\infty(M)$. Euclidean duality from Proposition~\ref{prop:estructura-escala-bessel-euclidiana}, uniform continuity of multiplication by $a_i$, and Hölder's inequality for sequences, applied to \eqref{eq:accion-distribucional-coordenadas-cap13}, give
\begin{equation*}
\begin{split}
|v(u)|
&\leq
C
\left(
\displaystyle\sum_{i\in I}
\left\|
\widetilde{(h_i u)\circ\phi_i^{-1}}
\right\|_{H^{r,p}(\mathbb R^n)}^p
\right)^{\frac{1}{p}}\\
&\hspace{12mm}\cdot
\left(
\displaystyle\sum_{i\in I}
\left\|
\widetilde{(h_i v)\circ\phi_i^{-1}}
\right\|_{H^{-r,p'}(\mathbb R^n)}^{p'}
\right)^{\frac{1}{p'}}\\
&\leq
C'\|u\|_{H^{r,p}(M)}\|v\|_{H^{-r,p'}(M)}.
\end{split}
\end{equation*}
By the density already proved, $u\mapsto v(u)$ has a unique continuous extension to $H^{r,p}(M)$. Define $\langle v,u\rangle_{r,p;\mathbf{g}}$ as the value of this extension. Then
\begin{equation}
\label{eq:extension-coincide-accion-distribucional-cap13}
\langle v,u\rangle_{r,p;\mathbf{g}}=v(u),
\qquad
v\in H^{-r,p'}(M),\quad u\in C_c^\infty(M),
\end{equation}
and, for two smooth functions, the right-hand side is $\displaystyle \int_Muv\,d\lambda_{\mathbf{g}}$. This explicitly shows how the measure induces the pairing and proves continuity of $\mathcal J_{r,p;\mathbf{g}}$.

It remains to prove surjectivity. Let $\Lambda\in(H^{r,p}(M))'$. By \eqref{eq:sintesis-continua-cap13}, $\Lambda\circ\mathcal R$ is continuous on $\ell^p(I;H^{r,p}(\mathbb R^n))$. Euclidean duality and Lemma~\ref{lem:dualidad-sumas-lq-banach} provide \[
w=(w_i)_{i\in I}
\in
\ell^{p'}(I;H^{-r,p'}(\mathbb R^n))
\] such that
\begin{equation}
\label{eq:representacion-funcional-sucesiones-cap13}
(\Lambda\circ\mathcal R)(z)
=
\displaystyle\sum_{i\in I}\langle w_i,z_i\rangle_{\mathbb R^n}
\end{equation}
for every $z=(z_i)_{i\in I}$, and
\[
\|w\|_{\ell^{p'}(I;H^{-r,p'}(\mathbb R^n))}
\leq
C\|\Lambda\|.
\]

Since $\operatorname{supp}(\chi_i)\Subset U_i$, the function \[
(\chi_i\circ\phi_i^{-1})a_i^{-1}
\] has a smooth extension by zero to $\mathbb R^n$. Denote it by
\[
m_i
:=
\widetilde{(\chi_i\circ\phi_i^{-1})a_i^{-1}}
\in C_c^\infty(\mathbb R^n).
\]
Its derivatives are uniformly bounded. Define
\begin{equation}
\label{eq:sintesis-representante-dual-cap13}
v
:=
\displaystyle\sum_{i\in I}h_i\,\mathcal E_i(m_iw_i).
\end{equation}
The proof of \eqref{eq:sintesis-continua-cap13}, with $-r,p'$ in place of $r,p$ and with the multipliers $m_i$, shows that $v\in H^{-r,p'}(M)$ and
\begin{equation}
\label{eq:cota-representante-dual-cap13}
\|v\|_{H^{-r,p'}(M)}
\leq
C\|w\|_{\ell^{p'}(I;H^{-r,p'}(\mathbb R^n))}
\leq
C'\|\Lambda\|.
\end{equation}

If $u\in C_c^\infty(M)$, only finitely many terms of \eqref{eq:sintesis-representante-dual-cap13} act on $u$. The definition of $\mathcal E_i$ and the rule for multiplying a distribution by a smooth function give
\begin{equation*}
\begin{split}
\langle v,u\rangle_{r,p;\mathbf{g}}
&=
v(u)\\
&=
\displaystyle\sum_{i\in I}
\left\langle
m_iw_i,
\widetilde{
a_i\bigl((h_i u)\circ\phi_i^{-1}\bigr)}
\right\rangle_{\mathbb R^n}\\
&=
\displaystyle\sum_{i\in I}
\left\langle
w_i,
\widetilde{(\chi_i h_i u)\circ\phi_i^{-1}}
\right\rangle_{\mathbb R^n}\\
&=
\displaystyle\sum_{i\in I}
\left\langle
w_i,
\widetilde{(h_i u)\circ\phi_i^{-1}}
\right\rangle_{\mathbb R^n}\\
&=
(\Lambda\circ\mathcal R)(\mathcal Su)
=
\Lambda(u).
\end{split}
\end{equation*}
The third equality uses $\chi_i=1$ near $\operatorname{supp}(h_i)$, and the penultimate equality uses \eqref{eq:representacion-funcional-sucesiones-cap13} and $\mathcal R\mathcal S=I$. Density of $C_c^\infty(M)$ extends the equality to all of $H^{r,p}(M)$, so $\mathcal J_{r,p;\mathbf{g}}$ is surjective. If $\mathcal J_{r,p;\mathbf{g}}(v)=0$, then $v$ vanishes on $C_c^\infty(M)$ and is the zero distribution. Thus, it is also injective. The direct estimate provides $C_D$, and \eqref{eq:cota-representante-dual-cap13} provides $c_D$ in \eqref{eq:dualidad-H-localizado}.
\end{proof}

\begin{corollary}[Functional connection with Amann's localization]
\label{cor:puente-funcional-amann-sobolev}
Let $(M,\mathbf{g})$ be a Riemannian manifold without boundary and with bounded geometry, let $\mathfrak K$ be the atlas constructed in Theorem~\ref{teo:geometria-acotada-implica-regularidad-uniforme-amann} or an equivalent uniformly regular atlas, and let $\mathcal Q_{\mathfrak K}=(U_i,\phi_i,h_i)_{i\in I}$ be the quadratic system of Proposition~\ref{prop:sistema-localizacion-amann}. For $s\in\mathbb R$ and $1<p<\infty$, the analysis and synthesis operators defined by \[
 \mathcal S_{\mathfrak K}u
 :=\left(\widetilde{(h_i u)\circ\phi_i^{-1}}\right)_{i\in I},
 \qquad
 \mathcal R_{\mathfrak K}v
 :=\displaystyle\sum_{i\in I}h_i\,\mathcal E_i(\widetilde\chi_i v_i)
\] are continuous, satisfy $\mathcal R_{\mathfrak K}\mathcal S_{\mathfrak K}
=\operatorname{id}_{H^{s,p}(M)}$, and there exist constants $c_{s,p},C_{s,p},D_{s,p}>0$ such that
\begin{equation}
\label{eq:comparacion-amann-localizacion-sobolev}
 c_{s,p}\|u\|_{H^{s,p}(M)}
 \leq
 \|\mathcal S_{\mathfrak K}u\|_{
   \ell^p(I;H^{s,p}(\mathbb R^m))}
 \leq
 C_{s,p}\|u\|_{H^{s,p}(M)}
\end{equation}
and
\[
 \|\mathcal R_{\mathfrak K}v\|_{H^{s,p}(M)}
 \leq D_{s,p}
 \|v\|_{\ell^p(I;H^{s,p}(\mathbb R^m))}.
\]
The constants depend on $m,s,p$, the multiplicity, and bounds on coordinate transitions and cutoffs through the order required by the local multipliers; they do not depend on $i$, $u$, or $v$.

For $1\leq q\leq\infty$, the localized norm with factors $h_i$ and the norm of $F^s_{p,q}(M)$ defined using the ordinary partition $\rho_i=h_i^2$ also satisfy both inequalities of Lemma~\ref{lem:equivalencia-localizacion-cuadratica}, with constants $c_F,C_F>0$ independent of $u$; when $q=p$, the same assertion holds for $B^s_{p,p}(M)$.
\end{corollary}

\begin{proof}
Proposition~\ref{prop:sistema-localizacion-amann} proves that $\mathcal Q_{\mathfrak K}$ is admissible. Lemma~\ref{lem:equivalencia-localizacion-cuadratica} provides the two inequalities in \eqref{eq:comparacion-amann-localizacion-sobolev} and the Triebel--Lizorkin and Besov comparisons. The distributional construction of $\mathcal E_i$, local finiteness of the sum, the synthesis bound, and the identity $\mathcal R\mathcal S=I$ were proved in Proposition~\ref{prop:retraccion-interpolacion-dualidad-geometria-acotada}; its proof uses only admissibility of the quadratic system and therefore applies unchanged to $\mathcal Q_{\mathfrak K}$. The stated dependencies are those of the multipliers, chart transitions, and the sum of at most $L$ terms appearing in that proof.
\end{proof}

Amann's uniform regularity thus describes the same norms through normalized charts. Comparisons with the quadratic system allow Euclidean estimates to be used with constants independent of the chart. To do the same for sections, we must also control the bundle connection and the transition matrices of its frames.

\begin{definition}[Bounded geometry of a vector bundle]
\label{def:haz-vectorial-geometria-acotada}
\index{vector bundle!of bounded geometry}
Let $(M,\mathbf{g})$ be a Riemannian manifold of bounded geometry without boundary, and let $\mathbf{E}\to M$ be a smooth vector bundle of finite rank equipped with a bundle metric $\mathbf{h}_{\mathbf{E}}$, Hermitian when the bundle is complex, and a compatible connection $\nabla^{\mathbf{E}}$. Denote the curvature of $\nabla^{\mathbf{E}}$ by \[
\mathbf{R}^{\mathbf{E}}\in\Gamma\bigl(\Lambda^2T^*M\otimes\operatorname{End}(\mathbf{E})\bigr)
\]. We say that $(\mathbf{E},\mathbf{h}_{\mathbf{E}},\nabla^{\mathbf{E}})$ is a \textit{bundle of bounded geometry} if, for each $k\in\mathbb N_0$, there exists $C_k^{\mathbf{E}}>0$ such that
\begin{equation}
\label{eq:curvatura-haz-acotada}
\displaystyle\sup_{x\in M}
\left|
(\nabla^{\Lambda^2T^*M\otimes\operatorname{End}(\mathbf{E})})^k\mathbf{R}^{\mathbf{E}}(x)
\right|_{\mathbf{g},\mathbf{h}_{\mathbf{E}}}
\leq C_k^{\mathbf{E}}.
\end{equation}
Here we use the product connection induced by the Levi--Civita connection and $\nabla^{\mathbf{E}}$.
For a fixed order $k\in\mathbb N_0$, if the base has bounded
geometry of order $k$, the bundle with its metric and connection
has \emph{bounded geometry of order $k$} when
\[
 \sup_{x\in M}
 \left|(\nabla^{\Lambda^2T^*M\otimes\operatorname{End}(\mathbf E)})^j
                   \mathbf R^{\mathbf E}(x)\right|_{\mathbf g,\mathbf h_{\mathbf E}}
 \leq C_j^{\mathbf E},\qquad 0\leq j\leq k,
\]
for constants $C_0^{\mathbf E},\ldots,C_k^{\mathbf E}>0$.
The connection being bounded is part of the bundle data.
\end{definition}

\begin{proposition}[Covariant derivatives and localization at integer orders]
\label{prop:covariantes-localizacion-orden-entero}
Let $(M,\mathbf{g})$ be a Riemannian manifold without boundary and with bounded geometry, let $\mathbf{E}\to M$ be a bundle of bounded geometry of rank $r_{\mathbf{E}}$, and let \[
\mathcal T^{\operatorname{geo}}
=(U_i,\phi_i,h_i)_{i\in I}
\] be a geodesic trivialization. On $U_i$, use the synchronous frame $\mathbf{e}_i=(\mathbf{e}_{i,1},\dots,\mathbf{e}_{i,r_{\mathbf{E}}})$ of Proposition~\ref{prop:marcos-sincronos-geodesicos-uniformes} and write
\[
\mathbf{u}=\displaystyle\sum_{a=1}^{r_{\mathbf{E}}}u_i^a\mathbf{e}_{i,a}.
\]
Set $W^{0,p}(M,\mathbf{E}):=L^p(M,\mathbf{E})$. If $m\in\mathbb N_0$ and $1\leq p<\infty$, there exist constants $c_{m,p},C_{m,p}>0$ such that, for every $\mathbf{u}\in W^{m,p}(M,\mathbf{E})$,
\begin{equation}
\label{eq:equivalencia-covariante-localizada-haz}
\begin{split}
c_{m,p}\|\mathbf{u}\|_{W^{m,p}(M,\mathbf{E})}
&\leq
\left(
\displaystyle\sum_{i\in I}\displaystyle\sum_{a=1}^{r_{\mathbf{E}}}
\left\|
\widetilde{(h_i u_i^a)\circ\phi_i^{-1}}
\right\|_{W^{m,p}(\mathbb R^n)}^p
\right)^{\frac{1}{p}}\\
&\leq
C_{m,p}\|\mathbf{u}\|_{W^{m,p}(M,\mathbf{E})}.
\end{split}
\end{equation}
If also $1<p<\infty$, then for the trivial rank-one bundle we obtain
\begin{equation}
\label{eq:W-H-entero-geometria-acotada}
W^{m,p}(M)=H^{m,p}(M)
\end{equation}
and there exist constants $c_{m,p}^{WH},C_{m,p}^{WH}>0$, depending on the
uniform geometric data, such that, for every function
$u\in W^{m,p}(M)=H^{m,p}(M)$,
\[
 c_{m,p}^{WH}\|u\|_{W^{m,p}(M)}
 \leq\|u\|_{H^{m,p}(M)}
 \leq C_{m,p}^{WH}\|u\|_{W^{m,p}(M)}.
\]
\end{proposition}

\begin{proof}
In the construction of Theorem~\ref{thm:trivializacion-geodesica}, we have $U_i=B_{\mathbf{g}}(p_i,2r)$. Since $3r<\operatorname{inj}(M)$ and $3r$ remains within the uniform radius of the metric criterion for bounded geometry, the normal chart extends to
\[
\widehat\phi_i\colon
\widehat U_i:=B_{\mathbf{g}}(p_i,3r)
\longrightarrow
B_{\mathrm{euc}}(0,3r),
\qquad
\widehat\phi_i\restriction_{U_i}=\phi_i.
\]
The uniform bounds on the metric, transitions between geodesic charts, and covariant derivatives of $\mathbf{R}^{\mathbf{E}}$ verify precisely the hypotheses of Proposition~\ref{prop:marcos-sincronos-geodesicos-uniformes}. In particular, the synchronous frame $\mathbf{e}_i$ is defined on a neighborhood of $\overline{U_i}$, and its connection coefficients have bounds independent of $i$. Moreover, each $(U_i,\phi_i)$ is a regular coordinate ball in the sense of Chapter~\ref{cap:sobolev-haces}. We may therefore apply Lemma~\ref{lema:meyers-serrin-haz-E} to every section $\mathbf{v}$ with compact support in $U_i$.

The constants in that lemma depend, through order $m$, on the components of the metric and its inverse, the volume density, the bundle metric in the chosen frame, and the connection coefficients. In the present normal charts and synchronous frames, all these quantities and their derivatives are uniformly bounded by the metric criterion for bounded geometry and Proposition~\ref{prop:marcos-sincronos-geodesicos-uniformes}. Since the synchronous frame is orthonormal,
\[
\left|\displaystyle\sum_{a=1}^{r_{\mathbf{E}}}z^a\mathbf{e}_{i,a}(x)\right|_{\mathbf{h}_{\mathbf{E}}(x)}=|z|
\qquad (x\in U_i,\ z\in\mathbb K^{r_{\mathbf{E}}}),
\]
and this identity allows the sum of the component norms in Lemma~\ref{lema:meyers-serrin-haz-E} to be replaced by their $\ell^p(I)$ sum. Thus, there exist $c_m^{\operatorname{loc}},C_m^{\operatorname{loc}}>0$, independent of $i$, such that
\begin{equation}
\label{eq:equivalencia-local-W-covariante-cap13}
\begin{split}
c_m^{\operatorname{loc}}
\left(
\displaystyle\sum_{a=1}^{r_{\mathbf{E}}}
\left\|
\widetilde{v_i^a\circ\phi_i^{-1}}
\right\|_{W^{m,p}(\mathbb R^n)}^p
\right)^{\frac{1}{p}}
&\leq
\left(
\displaystyle\sum_{q=0}^m
\|\nabla^q\mathbf{v}\|_{
L^p(U_i,T^{(0,q)}(TM)|_{U_i}\otimes \mathbf{E}|_{U_i})}^p
\right)^{\frac{1}{p}}\\
&\leq
C_m^{\operatorname{loc}}
\left(
\displaystyle\sum_{a=1}^{r_{\mathbf{E}}}
\left\|
\widetilde{v_i^a\circ\phi_i^{-1}}
\right\|_{W^{m,p}(\mathbb R^n)}^p
\right)^{\frac{1}{p}}.
\end{split}
\end{equation}

Apply this comparison to $\mathbf{v}=h_i \mathbf{u}$. For $q\in\{0,\dots,m\}$, the iterated Leibniz rule for the product connection is a finite sum, with coefficients depending only on $q$, of terms obtained by contracting
\[
\nabla^{q-\ell}h_i\otimes\nabla^\ell \mathbf{u},
\qquad \ell\in\{0,\dots,q\}.
\]
In particular,
\begin{equation}
\label{eq:leibniz-localizacion-covariante-cap13}
|\nabla^q(h_i \mathbf{u})|_{\mathbf{g},\mathbf{h}_{\mathbf{E}}}
\leq
C_q\displaystyle\sum_{\ell=0}^q
|\nabla^{q-\ell}h_i|_{\mathbf{g}}\,
|\nabla^\ell \mathbf{u}|_{\mathbf{g},\mathbf{h}_{\mathbf{E}}}.
\end{equation}
The derivatives of $h_i$ are uniformly bounded. Since each point belongs to at most $L$ supports, raising \eqref{eq:leibniz-localizacion-covariante-cap13} to the power $p$, summing over $i\in I$, and integrating gives
\begin{equation}
\label{eq:cota-localizacion-covariante-directa-cap13}
\displaystyle\sum_{i\in I}\displaystyle\sum_{q=0}^m
\|\nabla^q(h_i \mathbf{u})\|_{
L^p(U_i,T^{(0,q)}(TM)|_{U_i}\otimes \mathbf{E}|_{U_i})}^p
\leq
C\displaystyle\sum_{\ell=0}^m
\|\nabla^\ell \mathbf{u}\|_{L^p(M,T^{(0,\ell)}(TM)\otimes \mathbf{E})}^p.
\end{equation}
The first local inequality in \eqref{eq:equivalencia-local-W-covariante-cap13}, applied for each $i\in I$, proves the right-hand inequality of \eqref{eq:equivalencia-covariante-localizada-haz}.

For the opposite comparison, use the locally finite identity
\[
\mathbf{u}=\displaystyle\sum_{i\in I}h_i \mathbf{u}.
\]
For each $q\in\{0,\dots,m\}$, we have, distributionally and pointwise when $\mathbf{u}$ is smooth,
\[
\nabla^q \mathbf{u}
=
\displaystyle\sum_{i\in I}\nabla^q(h_i \mathbf{u}).
\]
At each point there are at most $L$ nonzero summands, so
\[
|\nabla^q \mathbf{u}|_{\mathbf{g},\mathbf{h}_{\mathbf{E}}}^p
\leq
L^{p-1}\displaystyle\sum_{i\in I}
|\nabla^q(h_i \mathbf{u})|_{\mathbf{g},\mathbf{h}_{\mathbf{E}}}^p.
\]
Integrate, sum over $q\in\{0,\dots,m\}$, and use the second inequality in \eqref{eq:equivalencia-local-W-covariante-cap13}. This gives the left-hand inequality of \eqref{eq:equivalencia-covariante-localizada-haz}.

To extend the comparisons to a section $\mathbf u\in W^{m,p}(M,\mathbf E)$, the Meyers--Serrin Theorem~\ref{meyers-serrin-haz} provides $\mathbf u_j\in\Gamma(\mathbf E)\cap W^{m,p}(M,\mathbf E)$ with $\mathbf u_j\to\mathbf u$ in $W^{m,p}(M,\mathbf E)$. Denote by $\mathcal S\mathbf v$ the family of components $\widetilde{(h_i v_i^a)\circ\phi_i^{-1}}$ appearing in \eqref{eq:equivalencia-covariante-localizada-haz}. The upper bound, applied to $\mathbf u_j-\mathbf u_k$, gives
\[
 \|\mathcal S\mathbf u_j-\mathcal S\mathbf u_k\|_
 {\ell^p(I\times\{1,\ldots,r_{\mathbf E}\};W^{m,p}(\mathbb R^n))}
 \leq C_{m,p}\|\mathbf u_j-\mathbf u_k\|_{W^{m,p}(M,\mathbf E)}.
\]
Completeness of this sequence space gives a limit $V$. For each fixed $i,a$, multiplication by $h_i$ and local equivalence of norms give convergence to $\widetilde{(h_i u_i^a)\circ\phi_i^{-1}}$ in $W^{m,p}(\mathbb R^n)$. Since projection onto a component is continuous, $V=\mathcal S\mathbf u$. Consequently, $\mathcal S\mathbf u_j\to\mathcal S\mathbf u$ in the norm of the entire family, not merely componentwise. Continuity of the norms allows limits to be taken in both inequalities applied to $\mathbf u_j$ and proves \eqref{eq:equivalencia-covariante-localizada-haz} for $\mathbf u$. The argument includes $p=1$: the factor $L^{p-1}$ is then one, and the sequence space remains complete. For the trivial rank-one bundle and $1<p<\infty$, Theorem~\ref{teo:identificaciones-H-F-W-entero} uniformly identifies the local norms $W^{m,p}(\mathbb R^n)$ and $H^{m,p}(\mathbb R^n)$ using its Euclidean constants. Combining them with $c_{m,p}$ and $C_{m,p}$ from \eqref{eq:equivalencia-covariante-localizada-haz} gives precisely $c_{m,p}^{WH}$ and $C_{m,p}^{WH}$ on $C_c^\infty(M)$. To justify compactly supported approximation as well, enumerate the partition and set $\displaystyle \chi_J:=\displaystyle\sum_{i=1}^J h_i$. Then $0\leq\chi_J\leq1$, each $\chi_J$ has compact support, and, on every fixed compact set, $\chi_J=1$ for sufficiently large $J$. Uniform local finiteness and the bounds on $\nabla^q h_i$ give
$\displaystyle \sup_{J\in\mathbb N}\|\nabla^q\chi_J\|_{L^\infty(M,T^{(0,q)}(TM))}<\infty$
for each $q$. For $0\leq q\leq m$, the Leibniz rule gives a constant independent of $J$ with
\[
 |\nabla_w^q((1-\chi_J)\mathbf u)|_{\mathbf g,\mathbf h_{\mathbf E}}
 \leq C_q\displaystyle\sum_{\ell=0}^q
 |\nabla^{q-\ell}(1-\chi_J)|_{\mathbf g}
 |\nabla_w^\ell\mathbf u|_{\mathbf g,\mathbf h_{\mathbf E}}.
\]
Every factor of the first type eventually vanishes on each fixed compact set and is uniformly bounded. The $p$th power of the right-hand side is dominated by $C\displaystyle\sum_{\ell=0}^m
|\nabla_w^\ell\mathbf u|_{\mathbf g,\mathbf h_{\mathbf E}}^p\in L^1(M)$. Theorem~\ref{convergencia dominada}, applied at each order and to their finite sum, proves $\chi_J\mathbf u\to\mathbf u$ in $W^{m,p}(M,\mathbf E)$. To retain compact support in the smooth approximation, fix $\psi_J\in C_c^\infty(M)$ equal to one near $\operatorname{supp}\chi_J$. If $\mathbf v_{J,k}\to\chi_J\mathbf u$ in $W^{m,p}(M,\mathbf E)$ is the Meyers--Serrin approximation, then
\[
 \|\psi_J\mathbf v_{J,k}-\chi_J\mathbf u\|_{W^{m,p}(M,\mathbf E)}
 \leq C_J\|\mathbf v_{J,k}-\chi_J\mathbf u\|_{W^{m,p}(M,\mathbf E)}.
\]
Choose $k(J)$ so that this last expression is less than $2^{-J}$. The sections $\psi_J\mathbf v_{J,k(J)}\in\Gamma_c(\mathbf E)$ converge to $\mathbf u$ in $W^{m,p}(M,\mathbf E)$, proving the asserted density. Density of $C_c^\infty(M)$ in $H^{m,p}(M)$, proved in Proposition~\ref{prop:retraccion-interpolacion-dualidad-geometria-acotada}, and the preceding approximation show that both spaces are completions of the same subspace with respect to these two explicit comparisons. The identification preserves their realizations in $\mathcal D'(M)$ and proves \eqref{eq:W-H-entero-geometria-acotada}.
\end{proof}

\begin{corollary}[Uniform control of the Bochner commutator]
\label{cor:control-uniforme-conmutador-bochner-geometria-acotada}
Let $(M,\mathbf{g})$ be a Riemannian manifold without boundary and with bounded geometry, and let $\mathbf{E}\to M$ be a bundle of bounded geometry. For $\ell\in\mathbb N_0$, let \[
\mathbf{E}_\ell:=T^{(0,\ell)}(TM)\otimes \mathbf{E}
\] carry the induced metric and product connection. Write
\[
\nabla_\ell:=\nabla^{\mathbf{E}_\ell}\colon
\Gamma(\mathbf{E}_\ell)\longrightarrow\Gamma(\mathbf{E}_{\ell+1}),
\qquad
\nabla^m:=\nabla_{\ell+m-1}\circ\cdots\circ\nabla_\ell,
\]
with the new index in the first covariant factor, and set $\Delta_{B,\ell}=(\nabla_\ell)_h^*\nabla_\ell
=-\operatorname{tr}_{\mathbf{g}}(\nabla_\ell^2)$, where the trace contracts the first two differentiation factors, for the positive-sign Bochner Laplacian. If $m\in\mathbb N$, then, for every $\mathbf{u}\in\Gamma_c(\mathbf{E}_\ell)$,
\begin{equation}
\label{eq:conmutador-bochner-uniforme-cap13}
\Delta_{B,\ell+m}(\nabla^m \mathbf{u})
-
\nabla^m(\Delta_{B,\ell}\mathbf{u})
=
\displaystyle\sum_{j=0}^m
\nabla^j\operatorname{Rm}*\nabla^{m-j}\mathbf{u}.
\end{equation}
Here $\operatorname{Rm}$ represents $\mathbf{R}^M$ or $\mathbf{R}^{\mathbf{E}}$ according to the term, and $*$ absorbs the tensor contractions and their universal coefficients. Moreover, there exists $C_{\ell,m}>0$ such that
\begin{equation}
\label{eq:cota-conmutador-bochner-uniforme-cap13}
\left|
\Delta_{B,\ell+m}(\nabla^m \mathbf{u})
-
\nabla^m(\Delta_{B,\ell}\mathbf{u})
\right|_{\mathbf{g},\mathbf{h}_{\mathbf{E}}}
\leq
C_{\ell,m}\displaystyle\sum_{q=0}^m
|\nabla^q \mathbf{u}|_{\mathbf{g},\mathbf{h}_{\mathbf{E}}}.
\end{equation}
\end{corollary}

\begin{proof}
Apply Lemma~\ref{lema conmutador laplaciano orden m} to the section $\mathbf{u}$ of the bundle $\mathbf{E}_\ell$. That result is formulated for the positive Bochner Laplacians of the induced tensor bundles and gives \eqref{eq:conmutador-bochner-uniforme-cap13} directly.

The curvature of $\mathbf{E}_\ell$ is the sum of the action of $\mathbf{R}^{\mathbf{E}}$ on the bundle factor and the action of $\mathbf{R}^M$ on each of the $\ell$ covariant factors. Bounded geometry of $M$ and $\mathbf{E}$ bounds all covariant derivatives of these curvatures. Every term on the right-hand side of \eqref{eq:conmutador-bochner-uniforme-cap13} is therefore bounded by a uniform multiple of $|\nabla^{m-j}\mathbf{u}|_{\mathbf{g},\mathbf{h}_{\mathbf{E}}}$. The finite sum gives \eqref{eq:cota-conmutador-bochner-uniforme-cap13}.
\end{proof}

\begin{proposition}[Perturbed chain of Bochner Laplacians]
\label{prop:cadena-perturbada-laplacianos-bochner-geometria-acotada}
Let $(M,\mathbf{g})$ be a Riemannian manifold without boundary and with bounded geometry, let $\mathbf{E}\to M$ be a bundle of bounded geometry, and let $N\in\mathbb N$. For $\ell\in\{0,\dots,N\}$, write \[
\mathbf{E}_\ell:=T^{(0,\ell)}(TM)\otimes \mathbf{E},
\qquad
 \nabla_\ell:=\nabla^{\mathbf{E}_\ell}\colon
\Gamma(\mathbf{E}_\ell)\longrightarrow\Gamma(\mathbf{E}_{\ell+1}),
\qquad
L_{\ell,0}:=(\nabla_\ell)_h^*\nabla_\ell
=-\operatorname{tr}_{\mathbf{g}}(\nabla_\ell^2)
\quad\text{on }\Gamma_c(\mathbf{E}_\ell),
\], and let $L_{\ell,p}$ be the positive generator of the covariant semigroup on $L^p(M,\mathbf{E}_\ell)$. There exist self-adjoint sections \[
\mathbf{V}_\ell\in\Gamma\bigl(\operatorname{End}(\mathbf{E}_\ell)\bigr),
\qquad \mathbf{V}_0=0,
\] and sections \[
\mathbf{W}_\ell\in
\Gamma\bigl(\operatorname{Hom}(\mathbf{E}_\ell,\mathbf{E}_{\ell+1})\bigr),
\qquad 0\leq\ell<N,
\] satisfying
\begin{equation}
\label{eq:conmutacion-cadena-bochner-perturbada-cap13}
\nabla_\ell(L_{\ell,0}+\mathbf{V}_\ell)\mathbf{u}
-(L_{\ell+1,0}+\mathbf{V}_{\ell+1})\nabla_\ell \mathbf{u}
=\mathbf{W}_\ell \mathbf{u},
\qquad \mathbf{u}\in\Gamma_c(\mathbf{E}_\ell).
\end{equation}
For each $a\in\mathbb N_0$, there exists a constant $C_{N,a}>0$ such that
\begin{equation}
\label{eq:cotas-coeficientes-cadena-bochner-cap13}
\displaystyle\sup_{\ell\in\{0,\ldots,N\}}
\displaystyle\sup_{x\in M}
|\nabla^a\mathbf{V}_\ell(x)|_{\mathbf{g},\mathbf{h}_{\mathbf{E}}}
+
\max_{\ell\in\{0,\ldots,N-1\}}
\displaystyle\sup_{x\in M}
|\nabla^a\mathbf{W}_\ell(x)|_{\mathbf{g},\mathbf{h}_{\mathbf{E}}}
\leq C_{N,a}.
\end{equation}

There exists $\mu_N>0$ such that, for $1\leq p<\infty$ and $\ell\in\{0,\dots,N\}$, the operator \begin{equation}
\label{eq:operador-B-cadena-bochner-cap13}
B_{\ell,p}
:=
\mu_N I+L_{\ell,p}+\mathbf{V}_\ell,
\qquad
\mathcal D(B_{\ell,p})=\mathcal D(L_{\ell,p}),
\end{equation} satisfies
\[
\|e^{-tB_{\ell,p}}\|_{\mathcal L(L^p(M,\mathbf{E}_\ell))}
\leq e^{-t},
\qquad t\geq0.
\]
In particular, $-(B_{\ell,p}-I)$ generates a contraction $C_0$-semigroup.
\end{proposition}

\begin{proof}
The bundle $\mathbf{E}_\ell$ is equipped with the product metric and connection. Its curvature is the sum of the action of $\mathbf{R}^{\mathbf{E}}$ on the factor $\mathbf{E}$ and the action of $\mathbf{R}^M$ on each of the $\ell$ covariant factors. Therefore, $\mathbf{E}_\ell$ has bounded geometry.

For a vector bundle $\mathbf F\to M$ equipped with a bundle metric $\mathbf h_{\mathbf F}$ and a compatible connection $\nabla^{\mathbf F}$, define the homomorphisms $\mathcal V_{\mathbf F}\in\Gamma(\operatorname{End}(T^*M\otimes\mathbf F))$ and $\mathcal W_{\mathbf F}\in\Gamma(\operatorname{Hom}(\mathbf F,T^*M\otimes\mathbf F))$ by
\begin{align}
 (\mathcal V_{\mathbf F}\mathbf a)(\mathbf X)
 &=\mathbf a(\operatorname{Ric}^{\sharp}\mathbf X)
   +2\displaystyle\sum_{i=1}^n\mathbf R^{\mathbf F}(\mathbf e_i,\mathbf X)
           \mathbf a(\mathbf e_i),
 \label{eq:potencial-curvatura-primer-orden-cap13}\\
 (\mathcal W_{\mathbf F}\mathbf v)(\mathbf X)
 &=\displaystyle\sum_{i=1}^n(\nabla_{\mathbf e_i}\mathbf R^{\mathbf F})
       (\mathbf e_i,\mathbf X)\mathbf v.
 \label{eq:termino-curvatura-orden-cero-cap13}
\end{align}
Here $\mathbf X\in T_xM$ is arbitrary, $\mathbf a\in T_x^*M\otimes\mathbf F_x$, $\mathbf v\in\mathbf F_x$, and $(\mathbf e_i)_{i=1}^n$ is an orthonormal basis of $T_xM$. The sums are metric contractions and hence independent of the basis. The exact first-order identity proved in Lemma~\ref{lema conmutador laplaciano orden m} states that
\begin{equation}
\label{eq:conmutador-primer-orden-bochner-cap13}
\nabla^{\mathbf{F}}\Delta_B^{\mathbf{F}}\mathbf{v}
-
\Delta_B^{T^*M\otimes \mathbf{F}}(\nabla^{\mathbf{F}}\mathbf{v})
=
\mathcal V_{\mathbf{F}}(\nabla^{\mathbf{F}}\mathbf{v})
+\mathcal W_{\mathbf{F}}\mathbf{v}.
\end{equation}
The homomorphism $\mathcal V_{\mathbf F}$ is self-adjoint. Indeed, $\operatorname{Ric}^{\sharp}$ is symmetric and, for the second term, the block taking component $i$ to component $j$ is $2\mathbf R^{\mathbf F}(\mathbf e_i,\mathbf e_j)$. Metric compatibility of the connection and antisymmetry of the curvature give
\[
 \bigl(\mathbf R^{\mathbf F}(\mathbf e_i,\mathbf e_j)\bigr)^*
 =-\mathbf R^{\mathbf F}(\mathbf e_i,\mathbf e_j)
 =\mathbf R^{\mathbf F}(\mathbf e_j,\mathbf e_i),
\]
which is exactly the relation between the adjoints of transposed blocks. If $\mathbf F$ has bounded geometry, all covariant derivatives of $\mathcal V_{\mathbf{F}}$ and $\mathcal W_{\mathbf{F}}$ are uniformly bounded, since both are obtained by contractions of the curvature of $\mathbf{F}$, of $\mathbf{R}^M$, and of its first covariant derivative.

Define recursively
\begin{equation}
\label{eq:definicion-potenciales-cadena-bochner-cap13}
\mathbf{V}_0:=0,
\qquad
\mathbf{V}_{\ell+1}
:=
I_{T^*M}\otimes \mathbf{V}_\ell+\mathcal V_{\mathbf{E}_\ell},
\qquad
\mathbf{W}_\ell:=\nabla \mathbf{V}_\ell+\mathcal W_{\mathbf{E}_\ell}.
\end{equation}
Here we use the canonical identification $\mathbf{E}_{\ell+1}=T^*M\otimes \mathbf{E}_\ell$. For $\mathbf{u}\in\Gamma_c(\mathbf{E}_\ell)$, the Leibniz rule gives
\[
\nabla_\ell(\mathbf{V}_\ell \mathbf{u})
=
(I_{T^*M}\otimes \mathbf{V}_\ell)\nabla_\ell \mathbf{u}
+(\nabla \mathbf{V}_\ell)\mathbf{u}.
\]
Adding this identity to \eqref{eq:conmutador-primer-orden-bochner-cap13}, with $\mathbf{F}=\mathbf{E}_\ell$, gives
\[
\begin{split}
&\nabla_\ell(L_{\ell,0}+\mathbf{V}_\ell)\mathbf{u}
-
(L_{\ell+1,0}+\mathbf{V}_{\ell+1})\nabla_\ell \mathbf{u}\\
&\qquad=
(\nabla \mathbf{V}_\ell+\mathcal W_{\mathbf{E}_\ell})\mathbf{u}
=
\mathbf{W}_\ell \mathbf{u}.
\end{split}
\]
This is \eqref{eq:conmutacion-cadena-bochner-perturbada-cap13}. Let us verify self-adjointness. The initial case is $\mathbf V_0=0$. If $\mathbf V_\ell$ is self-adjoint, then $I_{T^*M}\otimes\mathbf V_\ell$ is self-adjoint for the product metric of $\mathbf E_{\ell+1}=T^*M\otimes\mathbf E_\ell$; since $\mathcal V_{\mathbf E_\ell}$ is also self-adjoint, the recurrence \eqref{eq:definicion-potenciales-cadena-bochner-cap13} shows that $\mathbf V_{\ell+1}$ is self-adjoint. This completes the induction. To verify the bounds, fix $a\in\mathbb N_0$. Denote by $v_{\ell,a}$ a bound for $\|\nabla^a\mathbf V_\ell\|_\infty$, and by $A_{\ell,a},D_{\ell,a}$ the bounds for $\|\nabla^a\mathcal V_{\mathbf E_\ell}\|_\infty$ and $\|\nabla^a\mathcal W_{\mathbf E_\ell}\|_\infty$ already obtained from the curvatures. The connection preserves the identity tensor, so \eqref{eq:definicion-potenciales-cadena-bochner-cap13} implies
\[
 v_{0,a}=0,\qquad
 v_{\ell+1,a}\leq C_{n,a}\bigl(v_{\ell,a}+A_{\ell,a}\bigr),
 \qquad \ell\in\{0,\ldots,N-1\}.
\]
This recurrence starts at zero and has exactly $N$ steps. In particular, each $v_{\ell,a}$ is bounded by a finite sum of the $A_{j,a}$ with $j<\ell$, multiplied by powers of $C_{n,a}$. Moreover,
\[
 \|\nabla^a\mathbf W_\ell\|_\infty
 \leq C_{n,a}\bigl(v_{\ell,a+1}+D_{\ell,a}\bigr).
\]
Bounds at order $a+1$ follow from the same recurrence, independently of those at order $a$. Taking the maximum over the finitely many indices in the chain gives \eqref{eq:cotas-coeficientes-cadena-bochner-cap13}. This specifies which derivatives of the potentials enter each constant.

Theorem~\ref{teo:semigrupo-calor-covariante-Lp}, applied to $\mathbf{E}_\ell$, states that $-L_{\ell,p}$ generates a contraction semigroup. Multiplication by $\mathbf{V}_\ell$ is a bounded operator on $L^p(M,\mathbf{E}_\ell)$, and its norm is bounded by the supremum norm of its values as operators on the fibers. If $\mathbf{h}_{\mathbf{E}_\ell}$ denotes the product metric induced by $\mathbf{g}$ and $\mathbf{h}_{\mathbf{E}}$, set
\[
\|\mathbf V_\ell\|_{\infty,\mathrm{op}}
:=
\operatorname*{ess\,sup}_{x\in M}
\displaystyle\sup_{0\neq v\in(\mathbf{E}_\ell)_x}
\frac{|\mathbf{V}_\ell(x)v|_{\mathbf{h}_{\mathbf{E}_\ell}}}
{|v|_{\mathbf{h}_{\mathbf{E}_\ell}}}.
\]
Thus, the norm of the multiplication operator does not exceed
\[
\|\mathbf V_\ell\|_{\infty,\mathrm{op}}.
\]
Choose
\[
\mu_N
\geq
1+
\max_{\ell\in\{0,\ldots,N\}}
\|\mathbf V_\ell\|_{\infty,\mathrm{op}}.
\]
Proposition~\ref{prop:perturbacion-acotada-generador-semigrupo}, applied to the generator $-L_{\ell,p}$ and perturbation $-\mathbf{V}_\ell$, gives
\[
\|e^{-tB_{\ell,p}}\|_{\mathcal L(L^p(M,\mathbf{E}_\ell))}
\leq
\exp\left(
-t\bigl(\mu_N-
\|\mathbf V_\ell\|_{\infty,\mathrm{op}}\bigr)
\right)
\leq e^{-t}.
\]
Finally,
\[
e^{-t(B_{\ell,p}-I)}=e^t e^{-tB_{\ell,p}},
\]
so $\|e^{-t(B_{\ell,p}-I)}\|_{\mathcal L(L^p(M,\mathbf{E}_\ell))}\leq1$. This is the contraction semigroup generated by $-(B_{\ell,p}-I)$.
\end{proof}

The preceding differential identity will also be used for sections without compact support. For this, we need to identify the operator's $L^2(M,\mathbf F)$ domain, not only its local expression.

\begin{lemma}[Maximal domain of the perturbed Laplacian]
\label{lem:dominio-maximal-bochner-perturbado-uniforme}
Let $(\mathbf F,\mathbf h_{\mathbf F},\nabla^{\mathbf F})\to M$ be a bundle of bounded geometry over a manifold without boundary and with bounded geometry. Let $\mathbf V$ be a smooth, bounded, self-adjoint endomorphism, and let $B=L_{\mathbf F,2}+\mathbf V+\mu I$ be the realization through the energy form. If $\mathbf w,\mathbf f\in L^2(M,\mathbf F)$ and $B\mathbf w=\mathbf f$ distributionally, then $\mathbf w\in\mathcal D(B)$ and the equality holds in $L^2(M,\mathbf F)$.
\end{lemma}

\begin{proof}
Local elliptic regularity gives $\mathbf w\in H^2_{\mathrm{loc}}(M,\mathbf F)$. Take the cutoffs $\chi_J$ constructed in the proof of Proposition~\ref{prop:covariantes-localizacion-orden-entero}. They have compact support, satisfy $0\leq\chi_J\leq1$, converge locally to one, and
$|d\chi_J|\leq C$
with a constant independent of $J$. Test the equation against $\chi_J^2\mathbf w$, using a smooth approximation on its support. The real part of integration by parts gives
\[
\begin{aligned}
 \|\chi_J\nabla\mathbf w\|_{L^{2}(M,T^*M\otimes\mathbf F)}^2
 &\leq \|\mathbf f\|_{L^{2}(M,\mathbf F)}
          \|\mathbf w\|_{L^{2}(M,\mathbf F)}\\
 &\quad+(\|\mathbf V\|_\infty+|\mu|)
          \|\mathbf w\|_{L^{2}(M,\mathbf F)}^2\\
 &\quad+2\|\chi_J\nabla\mathbf w\|_{L^{2}(M,T^*M\otimes\mathbf F)}
          \,\bigl\||d\chi_J|\mathbf w\bigr\|_{L^{2}(M,\mathbf F)}.
\end{aligned}
\]
Young's inequality allows half of the energy term to be absorbed and gives a bound independent of $J$. Fatou implies $\nabla\mathbf w\in L^2(M,T^*M\otimes\mathbf F)$. Compactly supported density from the cited proposition then shows that $\mathbf w$ belongs to the form domain. For every test section $\boldsymbol\phi$,
\[
 \int_M\langle\nabla\mathbf w,\nabla\boldsymbol\phi\rangle
       \,d\lambda_{\mathbf g}
 +\int_M\langle(\mathbf V+\mu I)\mathbf w,\boldsymbol\phi\rangle
       \,d\lambda_{\mathbf g}
 =\int_M\langle\mathbf f,\boldsymbol\phi\rangle\,d\lambda_{\mathbf g}.
\]
All three terms are continuous for the form norm; density extends the equality to the entire form domain. The definition of the operator associated with a closed form gives $\mathbf w\in\mathcal D(B)$ and $B\mathbf w=\mathbf f$.
\end{proof}

The identification of fractional domains needed below will be obtained by interpolation. Before using functional calculus, it is useful to record interpolation of the covariant spaces at the two integer orders entering the square root.

\begin{proposition}[Interpolation between orders zero and two for bundles]
\label{prop:interpolacion-orden-cero-dos-haces-geometria-acotada}
Let $(M,\mathbf g)$ be a Riemannian manifold without boundary and with bounded geometry, and let $\mathbf F\to M$ be a finite-rank bundle of bounded geometry. If $1<p<\infty$, then
\begin{equation}
\label{eq:interpolacion-Lp-W2p-haz-cap13}
[\,L^p(M,\mathbf F),W^{2,p}(M,\mathbf F)\,]_{1/2}
=W^{1,p}(M,\mathbf F)
\end{equation}
with equivalence of norms.
\end{proposition}

\begin{proof}
If $\mathbf F$ is real, perform interpolation in the complexified bundle and restrict to the real subspace at the end; the connection, metric, and localization operators complexify without changing their norms on real sections. Use the same geodesic trivialization $\mathcal T^{\operatorname{geo}}=(U_i,\phi_i,h_i)_{i\in I}$ and synchronous frames $\mathbf e_i=(\mathbf e_{i,1},\ldots,\mathbf e_{i,r})$ as in Proposition~\ref{prop:covariantes-localizacion-orden-entero}. As in the definition of an admissible trivialization, choose functions $\chi_i\in C_c^\infty(U_i)$ such that
\[
 \chi_i=1\quad\text{on a neighborhood of }\operatorname{supp}(h_i),
\]
and whose derivatives in coordinates $\phi_i$ are uniformly bounded in $i$. Denote their extensions by zero by \[
 \widetilde\chi_i
 :=\widetilde{\chi_i\circ\phi_i^{-1}}\in C_c^\infty(\mathbb R^n)
\].

For $k\in\{0,1,2\}$, define the analysis operator
\[
 \mathcal S_{\mathbf F}\mathbf u
 :=
 \left(
 \widetilde{(h_i u_i^a)\circ\phi_i^{-1}}
 \right)_{(i,a)\in I\times\{1,\ldots,r\}},
 \qquad
 \mathbf u\restriction_{U_i}=\displaystyle\sum_{a=1}^r u_i^a\mathbf e_{i,a}.
\]
Proposition~\ref{prop:covariantes-localizacion-orden-entero} shows that $\mathcal S_{\mathbf F}$ is continuous from $W^{k,p}(M,\mathbf F)$ to
\[
 \ell^p\!\left(
 I\times\{1,\ldots,r\};W^{k,p}(\mathbb R^n)
 \right),
\]
with a norm depending on $k,p$ and the uniform data, but not on the chart.

Construct a synthesis operator that works simultaneously at all three orders. If $v=(v_i^a)_{i,a}$ is a family with finitely many nonzero components and each $v_i^a\in C_c^\infty(\mathbb R^n)$, set
\begin{equation}
\label{eq:sintesis-haz-interpolacion-012-cap13}
 \mathcal R_{\mathbf F}v
 :=
 \displaystyle\sum_{i\in I}
 \displaystyle\sum_{a=1}^r
 \left[
   \bigl((\widetilde\chi_i v_i^a)\restriction_{\phi_i(U_i)}\bigr)\circ\phi_i
 \right]\mathbf e_{i,a},
\end{equation}
where each summand is extended by zero outside $U_i$. The product with $\chi_i$ has compact support in $U_i$, so this extension is smooth. For a general family, the formula is first interpreted on each compact set: only uniformly finitely many charts contribute, and the sum is locally finite.

We prove continuity of synthesis. In a chart $U_j$, the components of the $i$th summand of \eqref{eq:sintesis-haz-interpolacion-012-cap13} are obtained from $v_i=(v_i^1,\ldots,v_i^r)$ through three operations: multiplication by a cutoff with uniformly bounded derivatives, a coordinate change $\phi_i\circ\phi_j^{-1}$, and multiplication by the transition matrix between the synchronous frames $\mathbf e_i$ and $\mathbf e_j$. The bounded-geometry bounds for synchronous frames and geodesic transitions imply that, for $k\in\{0,1,2\}$, these three operations are bounded on $W^{k,p}(\mathbb R^n,\mathbb K^r)$ with a constant independent of $i,j$. If $I_j$ denotes the set of indices whose charts meet $U_j$, there exists $L\geq1$ with $\#I_j\leq L$. By the inequality
\[
 \left(\displaystyle\sum_{i\in I_j}a_i\right)^p
 \leq L^{p-1}\displaystyle\sum_{i\in I_j}a_i^p,
 \qquad a_i\geq0,
\]
the local equivalence of Proposition~\ref{prop:covariantes-localizacion-orden-entero} gives
\begin{equation}
\label{eq:cota-sintesis-haz-interpolacion-012-cap13}
 \|\mathcal R_{\mathbf F}v\|_{W^{k,p}(M,\mathbf F)}
 \leq C_{k,p}
 \|v\|_{\ell^p(I\times\{1,\ldots,r\};W^{k,p}(\mathbb R^n))},
 \qquad k\in\{0,1,2\}.
\end{equation}
By density, $\mathcal R_{\mathbf F}$ extends continuously to all three sequence spaces.

If $v=\mathcal S_{\mathbf F}\mathbf u$, the equality $\chi_i=1$ near $\operatorname{supp}(h_i)$ makes the $i$th reconstructed term exactly $h_i\mathbf u$. Since $(h_i)_{i\in I}$ is a partition of unity,
\begin{equation}
\label{eq:retraccion-haz-interpolacion-012-cap13}
 \mathcal R_{\mathbf F}\mathcal S_{\mathbf F}\mathbf u
 =\displaystyle\sum_{i\in I}h_i\mathbf u=\mathbf u.
\end{equation}
Thus, the same analysis--synthesis pair is a coretraction--retraction at orders $0,1,2$.

By Corollary~\ref{teo:identificaciones-H-W-B-F} and interpolation of the Euclidean Bessel scale,
\[
 [L^p(\mathbb R^n),W^{2,p}(\mathbb R^n)]_{1/2}
 =W^{1,p}(\mathbb R^n)
\]
with equivalent norms. Proposition~\ref{prop:interpolacion-ellp-valores-banach} allows us to take the outer $\ell^p(I)$ sum and the finitely many bundle components. Now interpolate $\mathcal S_{\mathbf F}$ and $\mathcal R_{\mathbf F}$ between orders $0$ and $2$ and use \eqref{eq:retraccion-haz-interpolacion-012-cap13}. This gives both inclusions in \eqref{eq:interpolacion-Lp-W2p-haz-cap13}, and the endpoint operator norms provide the equivalence constants.
\end{proof}

The second tool is a global elliptic $L^p(M,\mathbf E_\ell)$ estimate with uniform constants. We prove it directly for the operators in the chain, making explicit the constants that will be used for the square root.

\begin{proposition}[$L^p(M,\mathbf E_\ell)$ domain of the chain operators]
\label{prop:dominio-Lp-cadena-bochner-geometria-acotada}
Under the hypotheses of Proposition~\ref{prop:cadena-perturbada-laplacianos-bochner-geometria-acotada}, let $1<p<\infty$. For each $\ell\in\{0,\ldots,N\}$,
\begin{equation}
\label{eq:dominio-Lp-B-cadena-W2p-cap13}
 \mathcal D(B_{\ell,p})=W^{2,p}(M,\mathbf E_\ell)
\end{equation}
and there exist constants $c_{\ell,p}^{(2)},C_{\ell,p}^{(2)}>0$ such that
\begin{equation}
\label{eq:equivalencia-grafica-B-W2p-cap13}
 c_{\ell,p}^{(2)}\|\mathbf u\|_{W^{2,p}(M,\mathbf E_\ell)}
 \leq
 \|B_{\ell,p}\mathbf u\|_{L^p(M,\mathbf E_\ell)}
 \leq
 C_{\ell,p}^{(2)}\|\mathbf u\|_{W^{2,p}(M,\mathbf E_\ell)}.
\end{equation}
The constants can be chosen uniformly as $\ell$ ranges over the finite set $\{0,\ldots,N\}$.
\end{proposition}

\begin{proof}
Fix $\ell$ and denote by \[
 \mathcal B_\ell:=\mu_NI+L_{\ell,0}+\mathbf V_\ell
\] the smooth differential expression whose realization is $B_{\ell,p}$. In a uniform geodesic chart and the synchronous frame of $\mathbf E_\ell$, this operator has the form
\begin{equation}
\label{eq:forma-local-B-cadena-Lp-cap13}
 \mathcal B_i
 =-g_i^{ab}(x)\partial_a\partial_b
  +A_i^a(x)\partial_a+C_i(x).
\end{equation}
The matrices $A_i^a,C_i$ and all derivatives we will need are uniformly bounded; the same holds for $g_i^{ab}$ and its derivatives. There are constants $0<\lambda_0\leq\Lambda_0<\infty$, independent of $i$, such that
\begin{equation}
\label{eq:elipticidad-uniforme-B-cadena-cap13}
 \lambda_0|\xi|^2
 \leq g_i^{ab}(x)\xi_a\xi_b
 \leq\Lambda_0|\xi|^2.
\end{equation}
These bounds follow from the metric criterion for bounded geometry and the bounds of Proposition \ref{prop:cadena-perturbada-laplacianos-bochner-geometria-acotada}.

We begin with a constant-coefficient estimate. For a point $x_0$ in the chart, set
\[
 q_{x_0}(\xi):=g_i^{ab}(x_0)\xi_a\xi_b,
 \qquad
 P_{x_0}:=-g_i^{ab}(x_0)\partial_a\partial_b.
\]
If $|\alpha|\leq2$, the symbol
\[
 m_{\alpha,x_0}(\xi)
 :=\frac{(i\xi)^\alpha}{1+q_{x_0}(\xi)}
\]
satisfies the hypotheses of Theorem~\ref{teo:multiplicadores-mikhlin} with constants independent of $x_0$ and $i$. Indeed, differentiating it $r$ times gives a finite sum of quotients whose denominators are powers of $1+q_{x_0}$ and whose numerators are polynomials in $\xi$ with coefficients bounded by $\lambda_0,\Lambda_0$; using \eqref{eq:elipticidad-uniforme-B-cadena-cap13} gives
\[
 |D_\xi^\beta m_{\alpha,x_0}(\xi)|
 \leq C_{\alpha,\beta}|\xi|^{-|\beta|},
 \qquad \xi\neq0,
\]
for every multi-index required by Mikhlin. Therefore,
\begin{equation}
\label{eq:estimacion-constante-W2p-cap13}
 \|v\|_{W^{2,p}(\mathbb R^n)}
 \leq C_p\|(I+P_{x_0})v\|_{L^p(\mathbb R^n)},
 \qquad v\in C_c^\infty(\mathbb R^n),
\end{equation}
with $C_p$ independent of $x_0$ and $i$.

We will also use an elementary intermediate-order inequality. For $0<\varepsilon\leq1$ and $j\in\{1,\ldots,n\}$, the symbol
\[
 \frac{i\varepsilon\xi_j}{1+\varepsilon^2|\xi|^2}
\]
and its Mikhlin derivatives are uniformly bounded in $\varepsilon$. Applying it to $(I-\varepsilon^2\Delta)v$ and summing over $j$ gives
\begin{equation}
\label{eq:interpolacion-derivada-primer-segundo-Lp-cap13}
 \|Dv\|_{L^p(\mathbb R^n)}
 \leq C_p\left(
 \varepsilon\|D^2v\|_{L^p(\mathbb R^n)}
 +\varepsilon^{-1}\|v\|_{L^p(\mathbb R^n)}
 \right),
 \qquad v\in W^{2,p}(\mathbb R^n).
\end{equation}

The first derivatives of $g_i^{ab}$ have a uniform bound. Choose $r>0$ independent of the chart and small enough that, on every Euclidean ball of radius $2r$ contained in a larger chart,
\[
 \|g_i^{ab}-g_i^{ab}(x_0)\|_{L^\infty(B_{\mathrm{euc}}(x_0,2r))}
 \leq\delta,
\]
where $\delta>0$ is fixed so that the corresponding term can be absorbed in \eqref{eq:estimacion-constante-W2p-cap13}. If $v\in C_c^\infty(B_{\mathrm{euc}}(x_0,2r))$, write
\[
 (I+P_{x_0})v
 =\mathcal B_i v
 +(g_i^{ab}-g_i^{ab}(x_0))\partial_a\partial_bv
 -A_i^a\partial_av+(I-C_i)v.
\]
The constant-coefficient estimate, the choice of $\delta$, and \eqref{eq:interpolacion-derivada-primer-segundo-Lp-cap13}, with $\varepsilon$ sufficiently small and then fixed, give
\begin{equation}
\label{eq:estimacion-local-B-W2p-cap13}
 \|v\|_{W^{2,p}(\mathbb R^n)}
 \leq C
 \left(
 \|\mathcal B_i v\|_{L^p(\mathbb R^n)}
 +\|v\|_{L^p(\mathbb R^n)}
 \right),
\end{equation}
with a constant independent of the chart.

The compactly supported estimate implies an interior estimate without assuming a global second-order norm a priori. Let $r\leq s<t\leq2r$, and choose $\zeta_{s,t}\in C_c^\infty(B_{\mathrm{euc}}(x_0,t))$ such that $\zeta_{s,t}=1$ on $B_{\mathrm{euc}}(x_0,s)$ and
\[
 |D\zeta_{s,t}|\leq C(t-s)^{-1},
 \qquad
 |D^2\zeta_{s,t}|\leq C(t-s)^{-2}.
\]
If $u$ is smooth on a neighborhood of $\overline{B_{\mathrm{euc}}(x_0,2r)}$, apply \eqref{eq:estimacion-local-B-W2p-cap13} to $\zeta_{s,t}u$. The Leibniz rule gives
\begin{align}
 \|D^2u\|_{L^p(B_s)}
 &\leq C\|\mathcal B_i u\|_{L^p(B_t)}
 +\frac{C}{t-s}\|Du\|_{L^p(B_t)}
 +C\left(1+\frac1{(t-s)^2}\right)\|u\|_{L^p(B_t)}.
\label{eq:estimacion-anillos-interior-B-W2p-cap13}
\end{align}
Here and in the next two formulas, we abbreviate $B_\rho=B_{\mathrm{euc}}(x_0,\rho)$.

We need to eliminate the first-order term on the right-hand side. The balls $B_t$, $r\leq t\leq2r$, are Lipschitz domains obtained from one another by dilations with factors in the compact interval $[1,2]$. Extend $u\restriction_{B_t}$ by Theorem~\ref{teo:extension-sobolev-uniforme-lipschitz} and apply \eqref{eq:interpolacion-derivada-primer-segundo-Lp-cap13} to the extension. The second-order estimate for the extension operator also contains $\|Du\|_{L^p(B_t)}$; multiply by the small parameter in the interpolation inequality and absorb that term into the left-hand side. Since the radii lie in the fixed interval $[r,2r]$, the extension and absorption constants are uniform. We thus obtain
\begin{equation}
\label{eq:interpolacion-local-bola-W2p-cap13}
 \|Du\|_{L^p(B_t)}
 \leq \eta\|D^2u\|_{L^p(B_t)}
 +C\eta^{-1}\|u\|_{L^p(B_t)},
 \qquad0<\eta\leq\eta_0.
\end{equation}
Take $\eta=\varepsilon(t-s)$, with $\varepsilon>0$ fixed and sufficiently small. If
\[
 \Phi(\rho):=\|D^2u\|_{L^p(B_\rho)},
\]
the preceding two estimates give numbers $0<\vartheta<1$ and $C>0$, independent of the chart, such that
\begin{equation}
\label{eq:iteracion-radios-W2p-cap13}
 \Phi(s)
 \leq\vartheta\Phi(t)
 +C\|\mathcal B_i u\|_{L^p(B_{2r})}
 +\frac{C}{(t-s)^2}\|u\|_{L^p(B_{2r})},
 \qquad r\leq s<t\leq2r.
\end{equation}
To iterate this inequality, choose $\rho\in(\sqrt\vartheta,1)$ and set
\[
 r_k:=2r-r\rho^k,
 \qquad k\in\mathbb N_0.
\]
Then $r_0=r$, the sequence $(r_k)_{k\in\mathbb N_0}$ is increasing,
$\displaystyle\lim_{k\to\infty}r_k=2r$, and $r_{k+1}-r_k=r(1-\rho)\rho^k$. Iterating \eqref{eq:iteracion-radios-W2p-cap13} $m$ times gives
\begin{align*}
 \Phi(r)
 &\leq \vartheta^m\Phi(r_m)
 +C\displaystyle\sum_{k=0}^{m-1}\vartheta^k
   \|\mathcal B_i u\|_{L^p(B_{2r})}\\
 &\quad+
 \frac{C}{r^2(1-\rho)^2}
 \displaystyle\sum_{k=0}^{m-1}\left(\frac{\vartheta}{\rho^2}\right)^k
 \|u\|_{L^p(B_{2r})}.
\end{align*}
Since $u$ is smooth on a neighborhood of $\overline{B_{2r}}$, $\Phi(2r)<\infty$, so $\vartheta^m\Phi(r_m)\to0$. Both geometric series converge, and, since $r$ is the same in all charts, we arrive at
\begin{equation}
\label{eq:estimacion-interior-B-W2p-cap13}
 \|u\|_{W^{2,p}(B_r)}
 \leq C\left(
 \|\mathcal B_i u\|_{L^p(B_{2r})}
 +\|u\|_{L^p(B_{2r})}
 \right).
\end{equation}
The first-order component is recovered from \eqref{eq:interpolacion-local-bola-W2p-cap13}; the order-zero terms are already on the right-hand side.

To obtain the global estimate, sum the interior estimates over a cover of bounded multiplicity. Use the geodesic trivialization of Theorem~\ref{thm:trivializacion-geodesica}: the balls $V_i:=B_{\mathbf g}(p_i,r)$ cover $M$, and the larger balls $U_i:=B_{\mathbf g}(p_i,2r)$ have uniformly bounded multiplicity. In normal coordinates and synchronous frames, uniform comparison of covariant and Euclidean norms, together with \eqref{eq:estimacion-interior-B-W2p-cap13}, gives, for every smooth section $\mathbf u$ with $\mathbf u,\mathcal B_\ell\mathbf u\in L^p(M,\mathbf E_\ell)$,
\begin{align*}
 \|\mathbf u\|_{W^{2,p}(M,\mathbf E_\ell)}^p
 &\leq C\displaystyle\sum_{i\in I}
 \|\mathbf u\|_{W^{2,p}(V_i,\mathbf E_\ell)}^p\\
 &\leq C\displaystyle\sum_{i\in I}\left(
 \|\mathcal B_\ell\mathbf u\|_{L^p(U_i,\mathbf E_\ell)}^p
 +\|\mathbf u\|_{L^p(U_i,\mathbf E_\ell)}^p
 \right)\\
 &\leq C\left(
 \|\mathcal B_\ell\mathbf u\|_{L^p(M,\mathbf E_\ell)}^p
 +\|\mathbf u\|_{L^p(M,\mathbf E_\ell)}^p
 \right).
\end{align*}
Consequently,
\begin{equation}
\label{eq:estimacion-global-grafica-B-W2p-cap13}
 \|\mathbf u\|_{W^{2,p}(M,\mathbf E_\ell)}
 \leq C\left(
 \|\mathcal B_\ell\mathbf u\|_{L^p(M,\mathbf E_\ell)}
 +\|\mathbf u\|_{L^p(M,\mathbf E_\ell)}
 \right).
\end{equation}
All sums have nonnegative terms; the last inequality uses only the uniform multiplicity bound of $(U_i)$. Thus, the argument is valid even when the global second-order norm is not known to be finite beforehand: the interior estimate makes each summand finite, and Tonelli allows them to be summed directly.

The realization $B_{\ell,p}$ has a bounded inverse given by
\[
 B_{\ell,p}^{-1}\mathbf f
 =\int_0^\infty e^{-tB_{\ell,p}}\mathbf f\,dt,
\]
and estimate $\|e^{-tB_{\ell,p}}\|\leq e^{-t}$ of Proposition~\ref{prop:cadena-perturbada-laplacianos-bochner-geometria-acotada} implies
\begin{equation}
\label{eq:inverso-B-cadena-norma-uno-cap13}
 \|B_{\ell,p}^{-1}\|_{\mathcal L(L^p(M,\mathbf E_\ell))}\leq1.
\end{equation}

Now let $\mathbf f\in L^p(M,\mathbf E_\ell)$. Choose $\mathbf f_j\in\Gamma_c(\mathbf E_\ell)$ with $\mathbf f_j\to\mathbf f$ in $L^p(M,\mathbf E_\ell)$ and set $\mathbf u_j=B_{\ell,p}^{-1}\mathbf f_j$. Consistency of the semigroups as $p$ varies shows that, since $\mathbf f_j\in L^2(M,\mathbf E_\ell)\cap L^p(M,\mathbf E_\ell)$,
\[
 \mathbf u_j=B_{\ell,2}^{-1}\mathbf f_j.
\]
Local elliptic regularity applied to $\mathcal B_\ell\mathbf u_j=\mathbf f_j$ shows that $\mathbf u_j$ is smooth. Applying \eqref{eq:estimacion-global-grafica-B-W2p-cap13} to $\mathbf u_j-\mathbf u_k$ and using \eqref{eq:inverso-B-cadena-norma-uno-cap13} gives
\[
 \|\mathbf u_j-\mathbf u_k\|_{W^{2,p}(M,\mathbf E_\ell)}
 \leq C\|\mathbf f_j-\mathbf f_k\|_{L^p(M,\mathbf E_\ell)}.
\]
By completeness, $\mathbf u_j$ converges in $W^{2,p}(M,\mathbf E_\ell)$ to a section $\mathbf v$. At the same time, $\mathbf u_j\to B_{\ell,p}^{-1}\mathbf f$ in $L^p(M,\mathbf E_\ell)$, so $\mathbf v=B_{\ell,p}^{-1}\mathbf f$. Thus, $\mathcal D(B_{\ell,p})\subseteq W^{2,p}(M,\mathbf E_\ell)$ and
\[
 \|\mathbf u\|_{W^{2,p}(M,\mathbf E_\ell)}
 \leq C\|B_{\ell,p}\mathbf u\|_{L^{p}(M,\mathbf E_\ell)},
 \qquad \mathbf u\in\mathcal D(B_{\ell,p}).
\]

For the opposite inclusion, let
$\mathbf u\in W^{2,p}(M,\mathbf E_\ell)$. Put
\[
 V_*:=\max_{0\leq\ell\leq N}\sup_{x\in M}
                  \|\mathbf V_\ell(x)\|_{\mathrm{op}}<\infty,
\]
where the operator norm is induced by the bundle metric.
The expression $\mathcal B_\ell=\mu_NI-\operatorname{tr}_{\mathbf g}\nabla^2
+\mathbf V_\ell$ and the Cauchy--Schwarz inequality for the trace give,
for every smooth section $\mathbf v$,
\[
 |\mathcal B_\ell\mathbf v|_{\mathbf h_{\mathbf E_\ell}}
 \leq\sqrt n\,|\nabla^2\mathbf v|_{\mathbf g,\mathbf h_{\mathbf E_\ell}}
       +(\mu_N+V_*)|\mathbf v|_{\mathbf h_{\mathbf E_\ell}}.
\]
The triangle inequality in $L^p(M,\mathbf E_\ell)$ yields
\begin{equation}
 \|\mathcal B_\ell\mathbf v\|_{L^p(M,\mathbf E_\ell)}
 \leq(\sqrt n+\mu_N+V_*)\|\mathbf v\|_{W^{2,p}(M,\mathbf E_\ell)}.
 \label{eq:cota-directa-B-W2p-cap13}
\end{equation}
Compactly supported density from
Proposition~\ref{prop:covariantes-localizacion-orden-entero} provides
$\mathbf u_j\in\Gamma_c(\mathbf E_\ell)$ with
$\mathbf u_j\to\mathbf u$ in $W^{2,p}(M,\mathbf E_\ell)$.
Estimate~\eqref{eq:cota-directa-B-W2p-cap13}, applied to
$\mathbf u_j-\mathbf u_k$, shows that
$(\mathcal B_\ell\mathbf u_j)$ is Cauchy in $L^p(M,\mathbf E_\ell)$.
Since $B_{\ell,p}$ is closed and agrees with $\mathcal B_\ell$ on
$\Gamma_c(\mathbf E_\ell)$, it follows that
$\mathbf u\in\mathcal D(B_{\ell,p})$. Passing to the limit in
\eqref{eq:cota-directa-B-W2p-cap13} gives the right-hand inequality in
\eqref{eq:equivalencia-grafica-B-W2p-cap13}. This proves
\eqref{eq:dominio-Lp-B-cadena-W2p-cap13} and equivalence of norms.
The constants are uniform in $\ell\in\{0,\ldots,N\}$ by the bounds in
Proposition~\ref{prop:cadena-perturbada-laplacianos-bochner-geometria-acotada}.
\end{proof}

To pass from order two to order one, we need a bound on imaginary powers. We will prove it within the book using a Bellman function and the Poisson semigroup. The algebraic part is separated into the following lemma.

\begin{lemma}[Bellman function for two Hilbert norms]
\label{lem:bellman-dos-normas-hilbert-cap13}
Let $2\leq p<\infty$ and $q=p/(p-1)$. There exist constants $\delta_p,c_p,C_p>0$ and a function $\mathcal B_p\colon H_1\times H_2\to[0,\infty)$, defined for any real Hilbert spaces $H_1,H_2$, with the following properties.
\begin{enumerate}[label=(\alph*)]
\item $\mathcal B_p$ is of class $C^1$, is $C^2$ away from the hypersurface $|\xi|^p=|\eta|^q$ and the axes, and
\[
 0\leq\mathcal B_p(\xi,\eta)
 \leq C_p\bigl(|\xi|^p+|\eta|^q\bigr).
\]
\item Its radial derivatives are nonnegative.
\item At every point of second-order differentiability, there exists $\tau=\tau(\xi,\eta)>0$ such that, for $a\in H_1$ and $b\in H_2$,
\begin{equation}
\label{eq:hessiana-bellman-pesos-cap13}
 d^2\mathcal B_p(\xi,\eta)[(a,b),(a,b)]
 \geq c_p\bigl(\tau|a|^2+\tau^{-1}|b|^2\bigr).
\end{equation}
In particular,
\begin{equation}
\label{eq:hessiana-bellman-producto-cap13}
 d^2\mathcal B_p(\xi,\eta)[(a,b),(a,b)]
 \geq2c_p|a|\,|b|.
\end{equation}
\item With the same choice of $\tau$,
\begin{equation}
\label{eq:gradiente-bellman-cap13}
 |d\mathcal B_p(\xi,\eta)[(a,b)]|^2
 \leq C_p\mathcal B_p(\xi,\eta)
 \bigl(\tau|a|^2+\tau^{-1}|b|^2\bigr).
\end{equation}
\end{enumerate}
The constants depend only on $p$, not on the dimensions of $H_1,H_2$.
\end{lemma}

\begin{proof}
For $u,v\geq0$, set
\[
 \beta_p(u,v)
 :=u^p+v^q+\delta_p\begin{cases}
 u^2v^{2-q},&u^p\leq v^q,\\[2mm]
 \displaystyle\frac2p u^p+
 \left(\frac2q-1\right)v^q,&u^p\geq v^q.
 \end{cases}
\]
On the boundary $u^p=v^q$, the two additional expressions agree. Indeed, there $v=u^{p-1}$ and $u^2v^{2-q}=u^p$, while $\frac2p+\frac2q-1=1$. The first derivatives also agree:
\[
 2uv^{2-q}=2u^{p-1},
 \qquad
 (2-q)u^2v^{1-q}=(2-q)v^{q-1}.
\]
Thus, $\beta_p$ is $C^1$. Define
\[
 \mathcal B_p(\xi,\eta):=\beta_p(|\xi|,|\eta|).
\]
The bound in item (a) follows from $u^2v^{2-q}\leq v^q$ when $u^p\leq v^q$ and from the definition in the other region. The radial derivatives are nonnegative if $0<\delta_p\leq1$.

We verify quantitative convexity. First suppose that $u^p\leq v^q$ and set
\[
 \tau:=v^{2-q}.
\]
The relation $u\leq v^{q-1}$ implies
\[
 u^{p-2}\leq v^{2-q},
 \qquad
 u^2v^{-q}\leq v^{q-2},
 \qquad
 uv^{1-q}\leq1.
\]
The second derivatives of the correction $u^2v^{2-q}$ are
\[
 2v^{2-q},\qquad
 2(2-q)uv^{1-q},\qquad
 -(2-q)(q-1)u^2v^{-q}.
\]
Consequently,
\[
 \partial_{uu}\beta_p\geq2\delta_p\tau,
\]
\[
 \partial_{vv}\beta_p
 \geq(q-1)\bigl(q-\delta_p(2-q)\bigr)\tau^{-1},
 \qquad
 |\partial_{uv}\beta_p|\leq2\delta_p(2-q).
\]
Choose $\delta_p>0$ small enough that $q-\delta_p(2-q)\geq q/2$ and the mixed term can be absorbed using
\[
 2|ab|\leq \varepsilon\tau a^2+
 \varepsilon^{-1}\tau^{-1}b^2.
\]
Thus, for some $c_p>0$,
\begin{equation}
\label{eq:hessiana-radial-bellman-region1-cap13}
 \partial_{uu}\beta_p\,a^2
 +2\partial_{uv}\beta_p\,ab
 +\partial_{vv}\beta_p\,b^2
 \geq c_p(\tau a^2+\tau^{-1}b^2).
\end{equation}
Moreover,
\[
 \frac{\partial_u\beta_p}{u}\geq2\delta_p\tau,
 \qquad
 \frac{\partial_v\beta_p}{v}\geq q\tau^{-1}.
\]

If $u^p\geq v^q$, set $\tau:=u^{p-2}$. The correction is diagonal, and the relation $v\leq u^{p-1}$, together with $q-2\leq0$, gives
\[
 v^{q-2}\geq u^{2-p}=\tau^{-1}.
\]
Hence,
\begin{equation}
\label{eq:hessiana-radial-bellman-region2-cap13}
 \partial_{uu}\beta_p\,a^2
 +\partial_{vv}\beta_p\,b^2
 \geq c_p(\tau a^2+\tau^{-1}b^2),
\end{equation}
and the same bounds, with another constant, hold for $(\partial_u\beta_p)/u$ and $(\partial_v\beta_p)/v$.

For a radial function of two Hilbert-space variables, the Hessian decomposes into its radial and orthogonal parts. Suppose $u=|\xi|>0$ and $v=|\eta|>0$, and define the signed radial components
\[
 a_{\mathrm r}:=\left\langle a,\frac{\xi}{u}\right\rangle_{H_1},
 \qquad
 b_{\mathrm r}:=\left\langle b,\frac{\eta}{v}\right\rangle_{H_2}.
\]
The corresponding orthogonal decompositions are
\[
 a=a_\parallel+a_\perp,
 \qquad a_\parallel=a_{\mathrm r}\frac{\xi}{u},
 \qquad \left\langle a_\perp,\xi\right\rangle_{H_1}=0,
\]
\[
 b=b_\parallel+b_\perp,
 \qquad b_\parallel=b_{\mathrm r}\frac{\eta}{v},
 \qquad \left\langle b_\perp,\eta\right\rangle_{H_2}=0.
\]
Differentiating the norm along the line $\xi+ta$ gives
\begin{align*}
 \left.\frac{d}{dt}|\xi+ta|\right|_{t=0}
 &=a_{\mathrm r},\\
 \left.\frac{d^2}{dt^2}|\xi+ta|\right|_{t=0}
 &=\frac{|a|^2-a_{\mathrm r}^2}{u}
 =\frac{|a_\perp|^2}{u}.
\end{align*}
The formulas for $\eta+tb$ are analogous. The chain rule then gives, in each region of second-order differentiability,
\begin{align*}
 d^2\mathcal B_p(\xi,\eta)[(a,b),(a,b)]
 &=(\partial_{uu}\beta_p)a_{\mathrm r}^2
   +2(\partial_{uv}\beta_p)a_{\mathrm r}b_{\mathrm r}
   +(\partial_{vv}\beta_p)b_{\mathrm r}^2\\
 &\quad+\frac{\partial_u\beta_p}{u}|a_\perp|^2
   +\frac{\partial_v\beta_p}{v}|b_\perp|^2,
\end{align*}
where the derivatives of $\beta_p$ are evaluated at $(u,v)$. Apply \eqref{eq:hessiana-radial-bellman-region1-cap13} to the real numbers $a_{\mathrm r},b_{\mathrm r}$ in the first region, and \eqref{eq:hessiana-radial-bellman-region2-cap13} in the second. Together with the bounds for $(\partial_u\beta_p)/u$ and $(\partial_v\beta_p)/v$, and decreasing $c_p$ if necessary, this gives
\begin{align*}
 d^2\mathcal B_p(\xi,\eta)[(a,b),(a,b)]
 &\geq c_p\bigl(
   \tau(a_{\mathrm r}^2+|a_\perp|^2)
   +\tau^{-1}(b_{\mathrm r}^2+|b_\perp|^2)\bigr)\\
 &=c_p\bigl(\tau|a|^2+\tau^{-1}|b|^2\bigr).
\end{align*}
This proves \eqref{eq:hessiana-bellman-pesos-cap13}. The inequality $\tau A^2+\tau^{-1}B^2\geq2AB$ gives \eqref{eq:hessiana-bellman-producto-cap13}.

The gradient bound remains. In the first region, the relations already used give
\[
 |\partial_u\beta_p|\leq C_p u\tau,
 \qquad
 |\partial_v\beta_p|\leq C_p v\tau^{-1},
\]
and therefore
\[
 \frac{|\partial_u\beta_p|^2}{\tau}
 +|\partial_v\beta_p|^2\tau
 \leq C_p(u^p+v^q)
 \leq C_p\beta_p(u,v).
\]
In the second region, write $v=u^{p-1}\rho$, with $0\leq\rho\leq1$. Since $q-1=1/(p-1)$, we have
\[
 v^{2q-2}u^{p-2}=u^p\rho^{2(q-1)}\leq u^p,
\]
and the explicit formulas for $\partial_u\beta_p$ and $\partial_v\beta_p$ give
\[
 |\partial_u\beta_p|^2\tau^{-1}
 +|\partial_v\beta_p|^2\tau
 \leq C_p(u^p+v^q).
\]
The Cauchy--Schwarz inequality with weights $\tau,\tau^{-1}$ proves \eqref{eq:gradiente-bellman-cap13}. Points on the axes follow by taking limits. When $p=2$, we have $q=2$, and the same function reduces to a piecewise quadratic form whose pieces agree globally; the preceding argument remains valid with $\tau=1$.
\end{proof}

\begin{theorem}[Laplace multipliers and imaginary powers of Bochner--Schrödinger operators]
\label{teo:potencias-imaginarias-bochner-schrodinger-autocontenido}
Let $(M,\mathbf g)$ be a complete Riemannian manifold without boundary, and let $\mathbf F\to M$ be a real or Hermitian vector bundle of finite rank equipped with a bundle metric and a compatible connection. Let $\mathbf W\in\Gamma(\operatorname{End}(\mathbf F))$ be smooth and self-adjoint, and suppose that
\begin{equation}
\label{eq:potencial-positivo-bochner-schrodinger-cap13}
 \langle\mathbf W(x)v,v\rangle\geq |v|^2
 \quad\text{for every }x\in M,
 \quad v\in\mathbf F_x.
\end{equation}
Let $B_2=L_{\mathbf F,2}+\mathbf W$ be the self-adjoint realization on $L^2(M,\mathbf F)$, and suppose that its semigroup extends consistently to semigroups $C_0$ on $L^r(M,\mathbf F)$, $1\leq r<\infty$, satisfying
\begin{equation}
\label{eq:decaimiento-calor-bochner-schrodinger-cap13}
 \|e^{-sB_r}\|_{\mathcal L(L^r(M,\mathbf F))}\leq e^{-s},
 \qquad s\geq0.
\end{equation}
If $\mathbf F$ is real, the complex powers and multipliers below are understood on the complexification, and the estimates are subsequently restricted to real sections whenever the corresponding operator preserves that subspace. Then, for every $1<p<\infty$, the spectral powers $B_2^{it}$ extend consistently to a $C_0$ group of bounded operators $B_p^{it}$ on $L^p(M,\mathbf F)$. There exist constants $C_p,\omega_p>0$, independent of the dimension of $M$ and the rank of $\mathbf F$, such that
\begin{equation}
\label{eq:potencias-imaginarias-bochner-schrodinger-cap13}
 \|B_p^{it}\|_{\mathcal L(L^p(M,\mathbf F))}
 \leq C_p e^{\omega_p|t|},
 \qquad t\in\mathbb R.
\end{equation}
More generally, if $\mathcal C_2:=B_2^{1/2}$ and, for $1<p<\infty$, $-\mathcal C_p$ denotes the generator of the subordinated Poisson semigroup on $L^p(M,\mathbf F)$, and if $\psi\in L^\infty(0,\infty)$, then the spectral multiplier
\begin{equation}
\label{eq:multiplicador-laplace-poisson-cap13}
 m_\psi(\mathcal C_2),
 \qquad
 m_\psi(\lambda)
 :=\lambda\int_0^\infty\psi(s)e^{-s\lambda}\,ds,
 \quad\lambda\geq1,
\end{equation}
extends to $L^p(M,\mathbf F)$ and
\begin{equation}
\label{eq:cota-multiplicador-laplace-poisson-cap13}
 \|m_\psi(\mathcal C_p)\|_{\mathcal L(L^p(M,\mathbf F))}
 \leq C_p\|\psi\|_{L^\infty(0,\infty)}.
\end{equation}
\end{theorem}

\begin{proof}
The proof is first carried out for compactly supported sections and then extended to $L^p(M,\mathbf F)$ by density. If the bundle is complex, regard it as a real bundle with inner product given by the real part of the Hermitian metric; the norms are unchanged.

\emph{The Poisson semigroup.} For $t>0$ and $1\leq r<\infty$, define
\begin{equation}
\label{eq:subordinacion-poisson-bochner-schrodinger-cap13}
 P_t^{(r)}\mathbf f
 :=\frac{t}{2\sqrt\pi}
 \int_0^\infty
 s^{-3/2}e^{-\frac{t^2}{4s}}e^{-sB_r}\mathbf f\,ds.
\end{equation}
The integral is a Bochner integral. The scalar subordination identity
\[
 e^{-t\sqrt\lambda}
 =\frac{t}{2\sqrt\pi}
 \int_0^\infty
 s^{-3/2}e^{-\frac{t^2}{4s}}e^{-s\lambda}\,ds,
 \qquad\lambda>0,
\]
follows by the substitution $u=t/(2\sqrt s)$ in the Gaussian integral. On $L^2(M,\mathbf F)$, the spectral theorem shows that $P_t^{(2)}=e^{-t\mathcal C_2}$. Moreover,
\eqref{eq:decaimiento-calor-bochner-schrodinger-cap13}
and the same formula with $\lambda=1$ give
\begin{equation}
\label{eq:decaimiento-poisson-bochner-schrodinger-cap13}
 \|P_t^{(r)}\|_{\mathcal L(L^r(M,\mathbf F))}\leq e^{-t}.
\end{equation}
Consistency of the heat semigroups implies consistency of the $P_t^{(r)}$. For $1<r<\infty$, denote the generator of $(P_t^{(r)})_{t\geq0}$ by $-\mathcal C_r$. On $L^2(M,\mathbf F)$, we have $\mathcal C_2=B_2^{1/2}$, and consistency shows that the resolvents and Laplace calculi of the operators $\mathcal C_r$ agree on intersections of the spaces $L^r(M,\mathbf F)$.

Let $\mathbf f\in\Gamma_c(\mathbf F)$. For $t>0$, $\mathbf U(t):=P_t^{(2)}\mathbf f$ belongs to the domain of every power of $B_2$, since $\lambda^me^{-t\sqrt\lambda}$ is bounded on $[1,\infty)$ for each $m$. Local elliptic regularity implies that $\mathbf U$ is smooth on $(0,\infty)\times M$, and spectral calculus gives
\begin{equation}
\label{eq:ecuacion-poisson-bochner-schrodinger-cap13}
 \partial_t^2\mathbf U
 =B_2\mathbf U
 =-\operatorname{tr}_{\mathbf g}(\nabla^2\mathbf U)+\mathbf W\mathbf U.
\end{equation}

\emph{Bilinear estimate.} Let $p,q$ be conjugate exponents. By symmetry, we may assume $p\geq2$; if $p<2$, interchange $p,q$ and the two sections. Let $\mathbf f\in\Gamma_c(\mathbf F)$ and $\mathbf g\in\Gamma_c(\mathbf F)$, and set
\[
 \mathbf U(t)=P_t\mathbf f,
 \qquad
 \mathbf V(t)=P_t\mathbf g.
\]
Apply Lemma~\ref{lem:bellman-dos-normas-hilbert-cap13} to the real fibers of $\mathbf F$ and write
\[
 b_\varepsilon(t,x)
 :=\mathcal B_p\bigl((\mathbf U(t,x),\varepsilon),
                    (\mathbf V(t,x),\varepsilon)\bigr)
   -\mathcal B_p((0,\varepsilon),(0,\varepsilon)),
 \qquad\varepsilon>0.
\]
Here we adjoin a parallel real direction to each fiber; thus, the two arguments of $\mathcal B_p$ never vanish. The function $\mathcal B_p$ is $C^1$ and piecewise $C^2$. On the set where the arguments have norm at least $\varepsilon$, its gradient is locally Lipschitz. Integrating by parts separately in the two regions $|\mathbf U_\varepsilon|^p<|\mathbf V_\varepsilon|^q$ and $|\mathbf U_\varepsilon|^p>|\mathbf V_\varepsilon|^q$, the interface terms cancel because the first derivatives of $\mathcal B_p$ agree there. The following chain rule therefore holds almost everywhere and distributionally. Equivalently, it can be obtained by approximating $\mathcal B_p$ by convolution on compact sets away from the axes and passing to the limit.

Let \[
 \mathscr L:=\partial_t^2+
 \operatorname{tr}_{\mathbf g}\nabla^2
\] be the scalar operator on $(0,\infty)\times M$. The Hessian formula for a radial function on a Hilbert fiber, together with \eqref{eq:ecuacion-poisson-bochner-schrodinger-cap13}, shows that the first-order terms in $\mathbf U,\mathbf V$ are given by the Hessian of $\mathcal B_p$ applied to the time and covariant derivatives. The order-zero terms are
\[
 \frac{\partial_u\beta_p}{|\mathbf U_\varepsilon|}
 \langle\mathbf W\mathbf U,\mathbf U\rangle
 +
 \frac{\partial_v\beta_p}{|\mathbf V_\varepsilon|}
 \langle\mathbf W\mathbf V,\mathbf V\rangle,
\]
and are nonnegative by \eqref{eq:potencial-positivo-bochner-schrodinger-cap13} and radial monotonicity of the Bellman function. Item (c) of the lemma, applied to each time derivative and each spatial orthonormal direction, provides a measurable function $\tau=\tau(t,x)>0$ such that
\begin{align}
 \mathscr L b_\varepsilon
 &\geq c_p\Bigl[
 \tau\bigl(|\partial_t\mathbf U|^2+|\nabla\mathbf U|^2\bigr)
 +\tau^{-1}\bigl(|\partial_t\mathbf V|^2+|\nabla\mathbf V|^2\bigr)
 \Bigr].
\label{eq:desigualdad-bellman-poisson-pesos-cap13}
\end{align}
In particular, by $\tau A^2+\tau^{-1}B^2\geq2AB$,
\begin{equation}
\label{eq:desigualdad-bellman-poisson-cap13}
 \mathscr L b_\varepsilon
 \geq2c_p
 \bigl(|\partial_t\mathbf U|^2+|\nabla\mathbf U|^2\bigr)^{1/2}
 \bigl(|\partial_t\mathbf V|^2+|\nabla\mathbf V|^2\bigr)^{1/2}
\end{equation}
almost everywhere.

Take the cutoff functions $\chi_j$ from Lemma~\ref{lem:funciones-corte-completitud-capitulo-calor}. Thus, $0\leq\chi_j\leq1$, $\chi_j\to1$ locally, and $\|d\chi_j\|_\infty\to0$. Integrate $t\chi_j^2\mathscr L b_\varepsilon$ over $M\times(a,R)$, with $0<a<R<\infty$. For fixed $j,a,R$, first let $\varepsilon\to0^{+}$. On this finite time cylinder and on the compact support of $d\chi_j$, the sections $\mathbf U,\mathbf V$ and their first derivatives are bounded. Since $\mathcal B_p$ is of class $C^1$, the function and its gradient are bounded on the compact set of arguments obtained by varying $0<\varepsilon\leq1$. This gives a dominating function independent of $\varepsilon$ for the boundary and cutoff terms. Dominated convergence therefore permits passage to the limit in every term. Only then let $a\to0^{+}$ and $R\to\infty$. To justify these time limits, observe that, for compactly supported data, spectral calculus on $L^2(M,\mathbf F)$ gives
\[
 \|\mathcal C_2^kP_t\mathbf f\|_{L^{2}(M,\mathbf F)}
 \leq\displaystyle\sup_{\lambda\geq1}\lambda^k e^{-t\lambda}\,\|\mathbf f\|_{L^{2}(M,\mathbf F)},
\]
and likewise for $\mathbf g$. In particular, terms multiplied by $t$ vanish at $t=\infty$; near zero, use the fact that smooth compactly supported sections belong to the domain of every power of $B_2$. Time integration produces the value at $t=0$. In the spatial variable, the divergence formula gives
\[
 \int_M\chi_j^2\operatorname{tr}_{\mathbf g}
       \nabla^2 b_\varepsilon\,d\lambda_{\mathbf g}
 =-2\int_M\chi_j\langle d\chi_j,db_\varepsilon\rangle_{\mathbf g}
       \,d\lambda_{\mathbf g}.
\]
The gradient bound \eqref{eq:gradiente-bellman-cap13}, the weighted inequality \eqref{eq:desigualdad-bellman-poisson-pesos-cap13}, and Young's inequality allow half of the weighted energy term to be absorbed, leaving
\begin{align}
&\int_0^\infty\!\int_M
 t\chi_j^2
 \bigl(|\partial_t\mathbf U|^2+|\nabla\mathbf U|^2\bigr)^{1/2}
 \bigl(|\partial_t\mathbf V|^2+|\nabla\mathbf V|^2\bigr)^{1/2}
 \,d\lambda_{\mathbf g}\,dt
\notag\\
&\qquad\leq C_p\int_M\chi_j^2
 \bigl(|\mathbf f|^p+|\mathbf g|^q\bigr)\,d\lambda_{\mathbf g}
\notag\\
&\qquad\quad+
 C_p\|d\chi_j\|_\infty^2
 \int_0^\infty t
 \left(
 \|\mathbf U(t)\|_{L^{p}(M,\mathbf F)}^p+\|\mathbf V(t)\|_{L^{q}(M,\mathbf F)}^q
 \right)dt.
\label{eq:bellman-con-cortes-cap13}
\end{align}
After the limit $\varepsilon\to0^{+}$ on the finite cylinder, the bound in item (a) of the lemma leaves only a multiple of $|\mathbf U|^p+|\mathbf V|^q$ in the cutoff term. Now use \eqref{eq:decaimiento-poisson-bochner-schrodinger-cap13}:
\[
 \|\mathbf U(t)\|_{L^{p}(M,\mathbf F)}^p\leq e^{-pt}\|\mathbf f\|_{L^{p}(M,\mathbf F)}^p,
 \qquad
 \|\mathbf V(t)\|_{L^{q}(M,\mathbf F)}^q\leq e^{-qt}\|\mathbf g\|_{L^{q}(M,\mathbf F)}^q.
\]
The time integral of the last term in \eqref{eq:bellman-con-cortes-cap13} is finite and independent of $j$. As $j\to\infty$, that term disappears and Fatou gives
\begin{equation}
\label{eq:estimacion-bilineal-poisson-cap13}
 \int_0^\infty\!\int_M
 t\,|\mathcal C_2P_t\mathbf f|\,|\mathcal C_2P_t\mathbf g|
 \,d\lambda_{\mathbf g}\,dt
 \leq C_p
 \left(\|\mathbf f\|_{L^{p}(M,\mathbf F)}^p+\|\mathbf g\|_{L^{q}(M,\mathbf F)}^q\right).
\end{equation}
We used only the time derivatives $\partial_tP_t=-\mathcal C_2P_t$ from the stronger estimate obtained above. Replace $\mathbf f$ by $a\mathbf f$ and $\mathbf g$ by $a^{-1}\mathbf g$. Choose $a>0$ so that $a^p\|\mathbf f\|_{L^{p}(M,\mathbf F)}^p=a^{-q}\|\mathbf g\|_{L^{q}(M,\mathbf F)}^q$. The left-hand side is unchanged, and we obtain the homogeneous form
\begin{equation}
\label{eq:estimacion-bilineal-poisson-homogenea-cap13}
 \int_0^\infty\!\int_M
 t\,|\mathcal C_2P_t\mathbf f|\,|\mathcal C_2P_t\mathbf g|
 \,d\lambda_{\mathbf g}\,dt
 \leq C_p\|\mathbf f\|_{L^{p}(M,\mathbf F)}\|\mathbf g\|_{L^{q}(M,\mathbf F)}.
\end{equation}

\emph{Laplace multipliers.} Let $\psi\in L^\infty(0,\infty)$. On $L^2(M,\mathbf F)$, define $m_\psi(\mathcal C_2)$ by spectral calculus. Since $\mathcal C_2\geq I$,
\[
 \mathbf g=\int_0^\infty \mathcal C_2P_r\mathbf g\,dr
\]
in $L^2(M,\mathbf F)$. For $\mathbf f,\mathbf g\in\Gamma_c(\mathbf F)$, Fubini, self-adjointness, and the semigroup law give
\begin{align*}
 \langle m_\psi(\mathcal C_2)\mathbf f,\mathbf g\rangle
 &=\int_0^\infty\!\int_0^\infty
 \psi(s)\langle \mathcal C_2P_s\mathbf f,\mathcal C_2P_r\mathbf g\rangle\,dr\,ds\\
 &=2\int_0^\infty
 \left(\int_0^{2t}\psi(s)\,ds\right)
 \langle \mathcal C_2P_t\mathbf f,\mathcal C_2P_t\mathbf g\rangle\,dt.
\end{align*}
The second equality uses the change of variable $t=(s+r)/2$, whose Jacobian is two. By \eqref{eq:estimacion-bilineal-poisson-homogenea-cap13},
\[
 |\langle m_\psi(\mathcal C_2)\mathbf f,\mathbf g\rangle|
 \leq4C_p\|\psi\|_\infty
 \|\mathbf f\|_{L^{p}(M,\mathbf F)}\|\mathbf g\|_{L^{q}(M,\mathbf F)}.
\]
Duality between $L^p(M,\mathbf F)$ and $L^q(M,\mathbf F)$ and density prove \eqref{eq:cota-multiplicador-laplace-poisson-cap13}.

\emph{Imaginary powers.} For $\tau\in\mathbb R$, set
\[
 \psi_\tau(s):=
 \frac{s^{-2i\tau}}{\Gamma(1-2i\tau)}.
\]
Since
\[
 \lambda\int_0^\infty s^{-2i\tau}e^{-s\lambda}\,ds
 =\Gamma(1-2i\tau)\lambda^{2i\tau},
\]
the multiplier of \eqref{eq:multiplicador-laplace-poisson-cap13} is $m_{\psi_\tau}(\mathcal C_2)=\mathcal C_2^{2i\tau}=B_2^{i\tau}$. Recall the elementary gamma-function bound we need. The beta--gamma identity and Euler's reflection formula give, for $y\in\mathbb R$,
\begin{equation}
\label{eq:modulo-gamma-imaginario-cap13}
 |\Gamma(1+iy)|^2
 =\Gamma(1+iy)\Gamma(1-iy)
 =\frac{\pi |y|}{\sinh(\pi|y|)},
\end{equation}
with value one at $y=0$ by continuity. The reflection formula used here follows from
\[
 \Gamma(z)\Gamma(1-z)
 =\int_0^\infty\frac{t^{z-1}}{1+t}\,dt
 =\frac{\pi}{\sin(\pi z)},
 \qquad 0<\operatorname{Re}z<1,
\]
where the first equality is the substitution $t=x/(1-x)$ in the beta--gamma identity, and the second follows by integrating $w^{z-1}/(1+w)$ over a keyhole contour around the positive axis; the two sides of the cut differ by the factor $e^{2\pi iz}$, and the only interior pole is $w=-1$. The identity, obtained initially for $0<\operatorname{Re}z<1$, extends meromorphically to values where both sides are defined. Apply it to $z=iy$ and then use $\Gamma(1+iy)=iy\Gamma(iy)$. From \eqref{eq:modulo-gamma-imaginario-cap13}, we deduce
\begin{equation}
\label{eq:cota-inversa-gamma-imaginario-cap13}
 \frac1{|\Gamma(1-2i\tau)|}
 \leq C(1+|\tau|)^{-1/2}e^{\pi|\tau|}
 \leq C e^{\pi|\tau|}.
\end{equation}
The Laplace multiplier bound then gives \eqref{eq:potencias-imaginarias-bochner-schrodinger-cap13}, for example with $\omega_p=\pi$ after modifying $C_p$.

Strong continuity on $\tau$ remains. On $L^2(M,\mathbf F)$, it follows immediately from spectral calculus. If $p\neq2$, choose $r\in(1,\infty)$ so that $p$ lies strictly between $2$ and $r$. For a compactly supported section, Hölder interpolates the $L^p(M,\mathbf F)$ norm between $L^2(M,\mathbf F)$ and $L^r(M,\mathbf F)$; the $L^2(M,\mathbf F)$ factor tends to zero as $\tau\to\tau_0$, and the $L^r(M,\mathbf F)$ factor remains bounded on compact intervals of $\tau$ by the estimate already proved. Density of $\Gamma_c(\mathbf F)$ in $L^p(M,\mathbf F)$ and local uniform boundedness of $B_p^{i\tau}$ extend continuity to all of $L^p(M,\mathbf F)$. The group law is verified on $L^2(M,\mathbf F)\cap L^p(M,\mathbf F)$ by spectral calculus and then on $L^p(M,\mathbf F)$ by density.
\end{proof}

The abstract step that converts the second-order information into a square root is recorded in the following lemma.

\begin{lemma}[Square root by interpolation and imaginary powers]
\label{lem:raiz-interpolacion-potencias-imaginarias-cap13}
Let $X$ be a reflexive complex Banach space, and let $-B$ be the generator of a $C_0$-semigroup satisfying $\|e^{-tB}\|\leq e^{-t}$ for $t\geq0$. Suppose that the negative powers $B^{-s}$, $s>0$, are defined by the Laplace calculus of that semigroup and that the imaginary powers form a $C_0$ group satisfying
\[
 \|B^{it}\|_{\mathcal L(X)}\leq Me^{\omega|t|}.
\]
Equip $\mathcal D(B)$ with the norm $\|Bx\|_X$ and $\mathcal D(B^{1/2})$ with $\|B^{1/2}x\|_X$. Then
\begin{equation}
\label{eq:interpolacion-dominio-raiz-cap13}
 [X,\mathcal D(B)]_{1/2}=\mathcal D(B^{1/2})
\end{equation}
with equivalence of norms. If $X$ is real, the same assertion holds after complexifying $X$ and $B$ and then restricting to the real subspace.
\end{lemma}

\begin{proof}
The identity is the case $a_0=0$, $a_1=1$ of domain interpolation. To prove it, fix an integer $N\geq2$. For $x\in\mathcal D(B^N)$ and $\varepsilon>0$, define on the strip $0\leq\operatorname{Re}z\leq1$
\[
 F_\varepsilon(z)
 :=e^{\varepsilon(z-1/2)^2}B^{1/2-z}x.
\]
Writing $B^{1/2-z}x=B^{1/2-z-N}B^Nx$, the power with negative real part is expressed by the Laplace integral; this proves holomorphy in the interior. On the two boundary lines,
\[
 F_\varepsilon(it)
 =e^{\varepsilon(it-1/2)^2}B^{-it}B^{1/2}x,
\]
\[
 B F_\varepsilon(1+it)
 =e^{\varepsilon(1/2+it)^2}B^{-it}B^{1/2}x.
\]
The bound on imaginary powers and the Gaussian factor show that $F_\varepsilon$ belongs to the Calderón space of the pair $(X,\mathcal D(B))$ and that
\[
 \|x\|_{[X,\mathcal D(B)]_{1/2}}
 \leq C\|B^{1/2}x\|_X.
\]
If $x\in\mathcal D(B^{1/2})$, set $x_\eta=(I+\eta B)^{-N}x$. Laplace functional calculus gives $x_\eta\in\mathcal D(B^N)$ and $B^{1/2}x_\eta\to B^{1/2}x$. The preceding estimate applied to differences proves that $x_\eta$ converges to $x$ in the interpolated space. Thus, $\mathcal D(B^{1/2})\subseteq[X,\mathcal D(B)]_{1/2}$.

For the reverse inclusion, let $F\in\mathscr F(X,\mathcal D(B))$ and $x=F(1/2)$. Fix $\eta>0$ and set $R_\eta=(I+\eta B)^{-N}$. The function
\[
 G_{\varepsilon,\eta}(z)
 :=e^{\varepsilon(z-1/2)^2}B^zR_\eta F(z)
\]
is holomorphic with values in $X$. On the boundary,
\[
 G_{\varepsilon,\eta}(it)
 =e^{\varepsilon(it-1/2)^2}B^{it}R_\eta F(it),
\]
and
\[
 G_{\varepsilon,\eta}(1+it)
 =e^{\varepsilon(1/2+it)^2}B^{it}R_\eta BF(1+it).
\]
The norms of $R_\eta$ remain uniformly bounded because the resolvent comes from the same Laplace calculus. The three-lines lemma gives, with a constant independent of $\eta$,
\[
 \|B^{1/2}R_\eta x\|_X
 \leq C\|F\|_{\mathscr F(X,\mathcal D(B))}.
\]
As $\eta\to0^{+}$, $R_\eta x\to x$ in $X$. By reflexivity, extract a subsequence along which $B^{1/2}R_\eta x$ converges weakly. The graph of the closed operator $B^{1/2}$ is a closed subspace of $X\times X$ and, by convexity, is also weakly closed. It follows that $x\in\mathcal D(B^{1/2})$ and
\[
 \|B^{1/2}x\|_X
 \leq C\|F\|_{\mathscr F(X,\mathcal D(B))}.
\]
Taking the infimum over $F$ proves the reverse inclusion and the second norm comparison.
\end{proof}

\begin{corollary}[Internal first-order estimate for the chain]
\label{cor:estimacion-primer-orden-cadena-bochner-interna}
Under the hypotheses of Proposition~\ref{prop:cadena-perturbada-laplacianos-bochner-geometria-acotada}, let $1<p<\infty$. For each $\ell\in\{0,\ldots,N\}$,
\begin{equation}
\label{eq:raiz-B-cadena-W1p-interna-cap13}
 \mathcal D(B_{\ell,p}^{1/2})=W^{1,p}(M,\mathbf E_\ell)
\end{equation}
with equivalence of norms. In particular, there exist $\widetilde a_{\ell,p},\widetilde b_{\ell,p}>0$ such that
\begin{equation}
\label{eq:cotas-raiz-B-cadena-W1p-interna-cap13}
\widetilde a_{\ell,p}
\left(
\|\mathbf z\|_{L^p(M,\mathbf E_\ell)}
+\|\nabla_\ell\mathbf z\|_{L^p(M,\mathbf E_{\ell+1})}
\right)
\leq
\|B_{\ell,p}^{1/2}\mathbf z\|_{L^p(M,\mathbf E_\ell)}
\leq
\widetilde b_{\ell,p}
\left(
\|\mathbf z\|_{L^p(M,\mathbf E_\ell)}
+\|\nabla_\ell\mathbf z\|_{L^p(M,\mathbf E_{\ell+1})}
\right).
\end{equation}
\end{corollary}

\begin{proof}
The choice of $\mu_N$ in Proposition~\ref{prop:cadena-perturbada-laplacianos-bochner-geometria-acotada} implies, pointwise,
\[
 \mu_N I+\mathbf V_\ell\geq I
\]
as a self-adjoint endomorphism of $\mathbf E_\ell$. The same proposition gives
\[
 \|e^{-tB_{\ell,r}}\|_{\mathcal L(L^r(M,\mathbf E_\ell))}\leq e^{-t}
 \qquad(1\leq r<\infty).
\]
Theorem~\ref{teo:potencias-imaginarias-bochner-schrodinger-autocontenido} therefore applies to $B_{\ell,p}$ and provides the imaginary powers required in Lemma~\ref{lem:raiz-interpolacion-potencias-imaginarias-cap13}. Proposition~\ref{prop:dominio-Lp-cadena-bochner-geometria-acotada} and Proposition~\ref{prop:interpolacion-orden-cero-dos-haces-geometria-acotada} give
\[
 \mathcal D(B_{\ell,p}^{1/2})
 =[L^p(M,\mathbf E_\ell),\mathcal D(B_{\ell,p})]_{1/2}
 =[L^p(M,\mathbf E_\ell),W^{2,p}(M,\mathbf E_\ell)]_{1/2}
 =W^{1,p}(M,\mathbf E_\ell).
\]
The norm comparisons in these three results, together with the definition of the $W^{1,p}(M,\mathbf E_\ell)$ norm, produce \eqref{eq:cotas-raiz-B-cadena-W1p-interna-cap13}.
\end{proof}

\begin{proposition}[Power estimates for the Bochner chain]
\label{prop:estimaciones-meyer-cadena-bochner-geometria-acotada}
Let $(M,\mathbf{g})$ be a Riemannian manifold without boundary and with bounded geometry, let $\mathbf{E}\to M$ be a bundle of bounded geometry, and let $N\in\mathbb N$ and $1<p<\infty$. Use the operators $B_{\ell,p}$ of Proposition~\ref{prop:cadena-perturbada-laplacianos-bochner-geometria-acotada}. For $k\in\mathbb N_0$ and $\ell\in\{0,\ldots,N\}$, define
\begin{equation}
\label{eq:dominio-bessel-cadena-bochner-cap13}
\mathscr D_{\ell,k}^p
:=
B_{\ell,p}^{-\frac{k}{2}}L^p(M,\mathbf{E}_\ell),
\qquad
\|\mathbf{u}\|_{\mathscr D_{\ell,k}^p}
:=
\|B_{\ell,p}^{\frac{k}{2}}\mathbf{u}\|_{L^p(M,\mathbf{E}_\ell)}.
\end{equation}
If also $k+\ell\leq N$, then
\begin{equation}
\label{eq:igualdad-dominios-cadena-bochner-cap13}
\mathscr D_{\ell,k}^p
=
W^{k,p}(M,\mathbf{E}_\ell),
\end{equation}
and there exist constants $c_{\ell,k,p},C_{\ell,k,p}>0$ such that
\begin{equation}
\label{eq:estimaciones-meyer-cadena-bochner-cap13}
\begin{split}
c_{\ell,k,p}
\displaystyle\sum_{j=0}^k
\|\nabla^j \mathbf{u}\|_{L^p(M,\mathbf{E}_{\ell+j})}
&\leq
\|B_{\ell,p}^{\frac{k}{2}}\mathbf{u}\|_{L^p(M,\mathbf{E}_\ell)}\\
&\leq
C_{\ell,k,p}
\displaystyle\sum_{j=0}^k
\|\nabla^j \mathbf{u}\|_{L^p(M,\mathbf{E}_{\ell+j})}.
\end{split}
\end{equation}
\end{proposition}

\begin{proof}
The negative power of \eqref{eq:dominio-bessel-cadena-bochner-cap13} is understood through
\begin{equation}
\label{eq:potencia-negativa-cadena-bochner-cap13}
B_{\ell,p}^{-\frac{s}{2}}\mathbf{f}
=
\frac{1}{\Gamma(\frac{s}{2})}
\int_0^\infty
t^{\frac{s}{2}-1}e^{-tB_{\ell,p}}\mathbf{f}\,dt,
\qquad s>0.
\end{equation}
Since $-(B_{\ell,p}-I)$ generates a contraction semigroup, Theorem~\ref{teo:escala-bessel-semigrupo-contraccion} guarantees that these powers are injective, their ranges are Banach spaces, and
\begin{equation}
\label{eq:monotonia-dominios-bessel-cadena-cap13}
\|B_{\ell,p}^{\frac{a}{2}}\mathbf{u}\|_{L^p(M,\mathbf{E}_\ell)}
\leq
\|B_{\ell,p}^{\frac{b}{2}}\mathbf{u}\|_{L^p(M,\mathbf{E}_\ell)},
\qquad 0\leq a\leq b,
\end{equation}
for $\mathbf{u}\in\mathscr D_{\ell,b}^p$.

The first-order estimate follows by combining elliptic regularity with domain interpolation. Proposition~\ref{prop:dominio-Lp-cadena-bochner-geometria-acotada} identifies the domain of $B_{\ell,p}$ with $W^{2,p}(M,\mathbf E_\ell)$, Theorem~\ref{teo:potencias-imaginarias-bochner-schrodinger-autocontenido} provides its imaginary powers, and Proposition~\ref{prop:interpolacion-orden-cero-dos-haces-geometria-acotada} identifies the interpolated space of order one. Corollary~\ref{cor:estimacion-primer-orden-cadena-bochner-interna} therefore gives
\begin{equation}
\label{eq:meyer-primer-orden-cadena-bochner-cap13}
\mathscr D_{\ell,1}^p=W^{1,p}(M,\mathbf{E}_\ell).
\end{equation}
Moreover, there exist constants $\widetilde a_{\ell,p},\widetilde b_{\ell,p}>0$ such that
\begin{equation}
\label{eq:cotas-meyer-primer-orden-cadena-bochner-cap13}
\widetilde a_{\ell,p}
\left(
\|\mathbf{z}\|_{L^p(M,\mathbf{E}_\ell)}
+\|\nabla_\ell \mathbf{z}\|_{L^p(M,\mathbf{E}_{\ell+1})}
\right)
\leq
\|B_{\ell,p}^{\frac{1}{2}}\mathbf{z}\|_{L^p(M,\mathbf{E}_\ell)}
\leq
\widetilde b_{\ell,p}
\left(
\|\mathbf{z}\|_{L^p(M,\mathbf{E}_\ell)}
+\|\nabla_\ell \mathbf{z}\|_{L^p(M,\mathbf{E}_{\ell+1})}
\right).
\end{equation}
Thus, the starting point of the induction is established by the uniform elliptic regularity, Bellman function, and complex interpolation proved above.

We now develop the commutation step needed at higher orders. Write
\[
S_\ell(t):=e^{-tB_{\ell,p}},
\qquad t\geq0.
\]
The covariant semigroups are consistent as $p$ varies. The Dyson--Phillips series of Proposition~\ref{prop:perturbacion-acotada-generador-semigrupo} shows that consistency is retained when the bounded potentials $\mathbf{V}_\ell$ and shift $\mu_NI$ are added.

Before commuting unbounded operators, fix a common regularizing core. Let $\varphi\in C_c^\infty((0,\infty))$ and define
\begin{equation}
\label{eq:regularizador-temporal-cadena-bochner-cap13}
\mathcal R_{\ell,\varphi}\mathbf{z}
:=
\int_0^\infty\varphi(t)S_\ell(t)\mathbf{z}\,dt.
\end{equation}
The integral is a Bochner integral. We need not assume that $S_\ell(t)\mathbf z$ is differentiable for arbitrary data in $L^p(M,\mathbf E_\ell)$. Extend $\varphi$ by zero to $\mathbb R$. For small $h>0$,
\[
 \frac{S_\ell(h)-I}{h}\mathcal R_{\ell,\varphi}\mathbf z
 =\int_0^\infty
   \frac{\varphi(t-h)-\varphi(t)}h S_\ell(t)\mathbf z\,dt
 \longrightarrow-\int_0^\infty\varphi'(t)S_\ell(t)\mathbf z\,dt.
\]
Convergence is in $L^p(M,\mathbf E_\ell)$, by dominated convergence on a compact interval. The definition of the generator $-B_{\ell,p}$ proves the formula for $m=1$; repeating it with $\varphi',\varphi'',\ldots$ gives
\begin{equation}
\label{eq:potencias-regularizador-temporal-cap13}
B_{\ell,p}^m\mathcal R_{\ell,\varphi}\mathbf{z}
=
\int_0^\infty\varphi^{(m)}(t)S_\ell(t)\mathbf{z}\,dt,
\qquad m\in\mathbb N_0.
\end{equation}
Therefore, $\displaystyle \mathcal R_{\ell,\varphi}\mathbf{z}\in\bigcap_{m\in\mathbb N}\mathcal D(B_{\ell,p}^m)$.

Choose a family $(\varphi_\varepsilon)_{\varepsilon>0}$ of nonnegative functions in $C_c^\infty((0,\varepsilon))$ with integral one. Strong continuity of the semigroup implies
\[
\mathcal R_{\ell,\varphi_\varepsilon}\mathbf{z}\longrightarrow \mathbf{z}
\quad\text{on }L^p(M,\mathbf{E}_\ell).
\]
This convergence holds in each space $\mathscr D_{\ell,\rho}^p$. Indeed, if $\mathbf{z}=B_{\ell,p}^{-\frac{\rho}{2}}\mathbf{f}$, Laplace functional calculus allows $B_{\ell,p}^{-\frac{\rho}{2}}$ to commute with $S_\ell(t)$, and
\[
\|\mathcal R_{\ell,\varphi_\varepsilon}\mathbf{z}-\mathbf{z}\|_{\mathscr D_{\ell,\rho}^p}
=
\|\mathcal R_{\ell,\varphi_\varepsilon}\mathbf{f}-\mathbf{f}\|_{L^p(M,\mathbf{E}_\ell)}
\longrightarrow0.
\]
As the common core, take
\[
 \mathscr C_{\ell,p}:=\bigcap_{m\in\mathbb N}\mathcal D(B_{\ell,p}^m).
\]
The preceding regularizations belong to this space and prove its density in every $\mathscr D_{\ell,\rho}^p$. Moreover, the power law shows that it is invariant under positive and negative powers of $B_{\ell,p}$ and under $S_\ell(t)$. This invariance will allow the same identity to be applied to $\mathbf z$ and $B_{\ell,p}^{a}\mathbf z$.

The geometric identities are first verified on regularizations of compactly supported data. Approximate $\mathbf{z}\in L^p(M,\mathbf{E}_\ell)$ by sections $\mathbf{z}_n\in\Gamma_c(\mathbf{E}_\ell)$. The sections $\mathcal R_{\ell,\varphi}\mathbf{z}_n$ belong simultaneously to $L^2(M,\mathbf E_\ell)$ and $L^p(M,\mathbf E_\ell)$ and, by \eqref{eq:potencias-regularizador-temporal-cap13}, to the domains of every power of $B_{\ell,2}$. Since $B_{\ell,2}$ has the same elliptic principal symbol as the Bochner Laplacian, local elliptic regularity shows that these sections are smooth. The differential identity \eqref{eq:conmutacion-cadena-bochner-perturbada-cap13} holds locally for them. It also justifies the realizations involved: if $\displaystyle \mathbf v\in\bigcap_{m\in\mathbb N}\mathcal D(B_{\ell,2}^m)$, the first-order estimate gives $\nabla_\ell\mathbf v,\nabla_\ell B_{\ell,2}\mathbf v\in L^2(M,\mathbf E_{\ell+1})$. By the differential identity,
\[
 B_{\ell+1,2}\nabla_\ell\mathbf v
 =\nabla_\ell B_{\ell,2}\mathbf v-\mathbf W_\ell\mathbf v\in L^2(M,\mathbf E_{\ell+1})
 \quad\hbox{in distributions}.
\]
Lemma~\ref{lem:dominio-maximal-bochner-perturbado-uniforme} implies $\nabla_\ell\mathbf v\in\mathcal D(B_{\ell+1,2})$. Applying the same argument to $S_\ell(t)\mathbf v$ and $B_{\ell,2}S_\ell(t)\mathbf v$ gives continuity in the graph norms needed to differentiate the product of semigroups. Formulas \eqref{eq:potencias-regularizador-temporal-cap13}, the first-order estimate \eqref{eq:cotas-meyer-primer-orden-cadena-bochner-cap13}, and continuity of $\mathbf{W}_\ell$ allow passage to the limit $n\to\infty$. Thus, every identity established below on $\mathcal R_{\ell,\varphi}\mathbf{z}_n$ first passes to $\mathcal R_{\ell,\varphi}\mathbf z$. Finally, let $\varphi=\varphi_\varepsilon$ and $\varepsilon\to 0^{+}$. Convergence in all graph norms of finite order extends the identities to $\mathscr C_{\ell,p}$.

Functional calculus for Bessel powers shows that, for $\rho\geq0$,
\begin{equation}
\label{eq:semigrupo-en-dominios-bessel-cadena-cap13}
\|S_\ell(t)\mathbf{z}\|_{\mathscr D_{\ell,\rho}^p}
\leq
e^{-t}\|\mathbf{z}\|_{\mathscr D_{\ell,\rho}^p}.
\end{equation}
Suppose that, for some $\rho\in\mathbb N_0$, multiplication by $\mathbf{W}_\ell$ is continuous from $\mathscr D_{\ell,\rho}^p$ to $\mathscr D_{\ell+1,\rho}^p$. For $t>0$, define
\begin{equation}
\label{eq:def-operador-duhamel-conmutador-cap13}
\mathcal K_\ell(t)\mathbf{z}
:=
\int_0^t
S_{\ell+1}(t-\tau)\mathbf{W}_\ell S_\ell(\tau)\mathbf{z}\,d\tau.
\end{equation}
The integral converges in $\mathscr D_{\ell+1,\rho}^p$ and
\begin{equation}
\label{eq:cota-duhamel-dominios-bessel-cap13}
\|\mathcal K_\ell(t)\mathbf{z}\|_{\mathscr D_{\ell+1,\rho}^p}
\leq
C_\rho t e^{-t}\|\mathbf{z}\|_{\mathscr D_{\ell,\rho}^p}.
\end{equation}

Let $\mathbf{z}\in\mathscr C_{\ell,p}$. On the smooth approximations described above, we may differentiate in $L^2(M,\mathbf E_{\ell+1})$ the function
\[
\tau\longmapsto
S_{\ell+1}(t-\tau)\nabla_\ell S_\ell(\tau)\mathbf{z},
\qquad 0<\tau<t.
\]
The identity
\[
\nabla_\ell B_{\ell,p}-B_{\ell+1,p}\nabla_\ell=\mathbf{W}_\ell,
\]
which is \eqref{eq:conmutacion-cadena-bochner-perturbada-cap13} with the same shift on both bundles, gives
\[
\frac{d}{d\tau}
\left[S_{\ell+1}(t-\tau)\nabla_\ell S_\ell(\tau)\mathbf{z}\right]
=
-S_{\ell+1}(t-\tau)\mathbf{W}_\ell S_\ell(\tau)\mathbf{z}.
\]
Integrate between $0$ and $t$ and pass to $L^p(M,\mathbf E_{\ell+1})$ by consistency. This gives
\begin{equation}
\label{eq:duhamel-conmutador-cadena-bochner-cap13}
S_{\ell+1}(t)\nabla_\ell \mathbf{z}-
\nabla_\ell S_\ell(t)\mathbf{z}
=
\mathcal K_\ell(t)\mathbf{z}.
\end{equation}
The first-order estimate and \eqref{eq:cota-duhamel-dominios-bessel-cap13} allow this identity to extend by density to $\mathscr D_{\ell,1}^p$; at higher orders it is interpreted in the $\mathscr D_{\ell+1,\rho}^p$ norm whenever both sides are defined. The expression on the right-hand side, however, is defined and continuous on all of $\mathscr D_{\ell,\rho}^p$.

Now consider the square root. For $\mathbf{z}\in\mathscr C_{\ell,p}$, the integral formula of Lemma~\ref{lem:formula-integral-raiz-semigrupo}, applied first with truncated limits, and \eqref{eq:duhamel-conmutador-cadena-bochner-cap13} give
\begin{align*}
&\frac{1}{2\sqrt\pi}
\int_\varepsilon^R t^{-\frac{3}{2}}
\bigl(\nabla_\ell \mathbf{z}-S_{\ell+1}(t)\nabla_\ell \mathbf{z}\bigr)\,dt\\
&\quad=
\nabla_\ell
\left[
\frac{1}{2\sqrt\pi}
\int_\varepsilon^R t^{-\frac{3}{2}}
\bigl(\mathbf{z}-S_\ell(t)\mathbf{z}\bigr)\,dt
\right]
-
\frac{1}{2\sqrt\pi}
\int_\varepsilon^R t^{-\frac{3}{2}}\mathcal K_\ell(t)\mathbf{z}\,dt.
\end{align*}
The first term on the right-hand side converges to $\nabla_\ell B_{\ell,p}^{\frac{1}{2}}\mathbf{z}$: the bracketed integral converges in $\mathscr D_{\ell,1}^p$ because $\mathbf{z}$ belongs to every domain, and we use \eqref{eq:cotas-meyer-primer-orden-cadena-bochner-cap13}. The second term converges in $\mathscr D_{\ell+1,\rho}^p$, since
\[
\int_0^\infty t^{-\frac{3}{2}}(t e^{-t})\,dt
=
\Gamma\!\left(\frac12\right)<\infty.
\]
Closedness of $B_{\ell+1,p}^{\frac{1}{2}}$ and the last assertion of Lemma~\ref{lem:formula-integral-raiz-semigrupo} show that $\nabla_\ell \mathbf{z}\in\mathscr D_{\ell+1,1}^p$ and that
\begin{equation}
\label{eq:defecto-raiz-cadena-bochner-cap13}
\Theta_{\ell,1}\mathbf{z}
:=
B_{\ell+1,p}^{\frac{1}{2}}\nabla_\ell \mathbf{z}
-
\nabla_\ell B_{\ell,p}^{\frac{1}{2}}\mathbf{z}
=
-\frac{1}{2\sqrt\pi}
\int_0^\infty t^{-\frac{3}{2}}\mathcal K_\ell(t)\mathbf{z}\,dt.
\end{equation}
By \eqref{eq:cota-duhamel-dominios-bessel-cap13},
\begin{equation}
\label{eq:cota-defecto-raiz-dominios-cap13}
\|\Theta_{\ell,1}\mathbf{z}\|_{\mathscr D_{\ell+1,\rho}^p}
\leq
C_\rho'\|\mathbf{z}\|_{\mathscr D_{\ell,\rho}^p}.
\end{equation}
Since $\mathscr C_{\ell,p}$ is dense in $\mathscr D_{\ell,\rho}^p$, the integral on the right defines the unique extension of $\Theta_{\ell,1}$ to that space.

For $j\in\mathbb N$, write formally
\begin{equation}
\label{eq:defecto-potencias-cadena-bochner-cap13}
\Theta_{\ell,j}\mathbf{z}
:=
B_{\ell+1,p}^{\frac{j}{2}}\nabla_\ell \mathbf{z}
-
\nabla_\ell B_{\ell,p}^{\frac{j}{2}}\mathbf{z}.
\end{equation}
Before using this expression, we must verify that its first term is defined. We do this simultaneously with the defect estimate.

Suppose that, for some $R\in\mathbb N_0$, multiplication by $\mathbf{W}_\ell$ is continuous from $\mathscr D_{\ell,\rho}^p$ to $\mathscr D_{\ell+1,\rho}^p$ for every $\rho\in\{0,\dots,R\}$. We claim that, if $j\geq1$ and $j-1\leq R$, the following properties hold.
\begin{enumerate}[label=(\roman*)]
\item For each $\mathbf{z}\in\mathscr C_{\ell,p}$, $\nabla_\ell \mathbf{z}\in\mathscr D_{\ell+1,j}^p$, and the identity \eqref{eq:defecto-potencias-cadena-bochner-cap13} is well defined.
\item If $\rho\in\mathbb N_0$ and $\rho+j-1\leq R$, then $\Theta_{\ell,j}$ extends uniquely to a continuous operator
\[
\Theta_{\ell,j}\colon
\mathscr D_{\ell,\rho+j-1}^p
\longrightarrow
\mathscr D_{\ell+1,\rho}^p
\]
and there exists $C_{\ell,j,\rho,p}>0$ such that
\begin{equation}
\label{eq:cota-defecto-potencias-refinada-cap13}
\|\Theta_{\ell,j}\mathbf{z}\|_{\mathscr D_{\ell+1,\rho}^p}
\leq
C_{\ell,j,\rho,p}
\|\mathbf{z}\|_{\mathscr D_{\ell,\rho+j-1}^p}.
\end{equation}
\end{enumerate}

For $j=1$, both assertions follow from \eqref{eq:defecto-raiz-cadena-bochner-cap13} and \eqref{eq:cota-defecto-raiz-dominios-cap13}. Moreover, the square-root identity extends to the mapping
\begin{equation}
\label{eq:gradiente-dos-a-uno-cadena-bochner-cap13}
\nabla_\ell\colon
\mathscr D_{\ell,2}^p
\longrightarrow
\mathscr D_{\ell+1,1}^p.
\end{equation}
Indeed, if $\mathbf{z}_\nu\in\mathscr C_{\ell,p}$ converges to $\mathbf{z}$ in $\mathscr D_{\ell,2}^p$, the first-order estimate shows that $\nabla_\ell \mathbf{z}_\nu\to\nabla_\ell \mathbf{z}$ in $L^p(M,\mathbf E_{\ell+1})$. On the other hand,
\[
B_{\ell+1,p}^{\frac{1}{2}}\nabla_\ell \mathbf{z}_\nu
=
\nabla_\ell B_{\ell,p}^{\frac{1}{2}}\mathbf{z}_\nu
+
\Theta_{\ell,1}\mathbf{z}_\nu
\]
is Cauchy in $L^p(M,\mathbf E_{\ell+1})$: the first term is controlled by applying \eqref{eq:cotas-meyer-primer-orden-cadena-bochner-cap13} to $B_{\ell,p}^{\frac{1}{2}}(\mathbf{z}_\nu-\mathbf{z}_\mu)$, and the second by \eqref{eq:cota-defecto-raiz-dominios-cap13} with $\rho=0$. Closedness of $B_{\ell+1,p}^{\frac{1}{2}}$ proves \eqref{eq:gradiente-dos-a-uno-cadena-bochner-cap13} and permits passage to the limit in the identity.

Proceed by induction on $j$. Suppose both assertions have been proved through order $j-1$, with $j\geq2$, and let $\mathbf{z}\in\mathscr C_{\ell,p}$. The induction hypothesis gives
\begin{equation}
\label{eq:identidad-defecto-orden-anterior-cap13}
B_{\ell+1,p}^{\frac{j-1}{2}}\nabla_\ell \mathbf{z}
=
\nabla_\ell B_{\ell,p}^{\frac{j-1}{2}}\mathbf{z}
+
\Theta_{\ell,j-1}\mathbf{z}.
\end{equation}
The first term on the right belongs to $\mathscr D_{\ell+1,1}^p$ by \eqref{eq:gradiente-dos-a-uno-cadena-bochner-cap13}, applied to $B_{\ell,p}^{\frac{j-1}{2}}\mathbf{z}$. The second also belongs to that space by \eqref{eq:cota-defecto-potencias-refinada-cap13} with order $j-1$ and $\rho=1$; here we use $j-1\leq R$. Thus, the left-hand side of \eqref{eq:identidad-defecto-orden-anterior-cap13} belongs to the domain of $B_{\ell+1,p}^{\frac{1}{2}}$. This proves that $\nabla_\ell \mathbf{z}\in\mathscr D_{\ell+1,j}^p$.

Now apply $B_{\ell+1,p}^{\frac{1}{2}}$ to \eqref{eq:identidad-defecto-orden-anterior-cap13} and use the square-root identity on $B_{\ell,p}^{\frac{j-1}{2}}\mathbf{z}$. We obtain the recursion
\begin{equation}
\label{eq:recursion-defectos-potencias-cap13}
\Theta_{\ell,j}\mathbf{z}
=
\Theta_{\ell,1}\bigl(B_{\ell,p}^{\frac{j-1}{2}}\mathbf{z}\bigr)
+
B_{\ell+1,p}^{\frac{1}{2}}\Theta_{\ell,j-1}\mathbf{z}.
\end{equation}
If $\rho+j-1\leq R$, estimate the first sum using \eqref{eq:cota-defecto-raiz-dominios-cap13} in the $\mathscr D_{\ell+1,\rho}^p$ norm, and the second using the induction hypothesis in the $\mathscr D_{\ell+1,\rho+1}^p$ norm. Thus,
\[
\begin{split}
\|\Theta_{\ell,j}\mathbf{z}\|_{\mathscr D_{\ell+1,\rho}^p}
&\leq
C_\rho
\|B_{\ell,p}^{\frac{j-1}{2}}\mathbf{z}\|_{\mathscr D_{\ell,\rho}^p}
+
\|\Theta_{\ell,j-1}\mathbf{z}\|_{\mathscr D_{\ell+1,\rho+1}^p}\\
&\leq
C_{\ell,j,\rho,p}
\|\mathbf{z}\|_{\mathscr D_{\ell,\rho+j-1}^p}.
\end{split}
\]
Density of $\mathscr C_{\ell,p}$ in the source space provides the unique extension and concludes the induction. In particular, taking $\rho=0$ gives
\begin{equation}
\label{eq:cota-defecto-potencias-cadena-bochner-cap13}
\|\Theta_{\ell,j}\mathbf{z}\|_{L^p(M,\mathbf{E}_{\ell+1})}
\leq
C_{\ell,j,p}
\|B_{\ell,p}^{\frac{j-1}{2}}\mathbf{z}\|_{L^p(M,\mathbf{E}_\ell)}.
\end{equation}

We now deduce the connection's mapping property between domains. Let $\mathbf{z}\in\mathscr D_{\ell,j+1}^p$, and choose $\mathbf{z}_\nu\in\mathscr C_{\ell,p}$ with $\mathbf{z}_\nu\to \mathbf{z}$ in $\mathscr D_{\ell,j+1}^p$. The first-order estimate gives $\nabla_\ell \mathbf{z}_\nu\to\nabla_\ell \mathbf{z}$ in $L^p(M,\mathbf E_{\ell+1})$. On the other hand,
\[
B_{\ell+1,p}^{\frac{j}{2}}\nabla_\ell \mathbf{z}_\nu
=
\nabla_\ell B_{\ell,p}^{\frac{j}{2}}\mathbf{z}_\nu
+
\Theta_{\ell,j}\mathbf{z}_\nu.
\]
The right-hand side is Cauchy in $L^p(M,\mathbf E_{\ell+1})$: apply the first-order estimate to $B_{\ell,p}^{\frac{j}{2}}(\mathbf{z}_\nu-\mathbf{z}_\mu)$ for the first term and use \eqref{eq:cota-defecto-potencias-cadena-bochner-cap13} for the second. Closedness of $B_{\ell+1,p}^{\frac{j}{2}}$ implies $\nabla_\ell \mathbf{z}\in\mathscr D_{\ell+1,j}^p$ and permits passage to the limit in the identity. Moreover, the first-order estimate and \eqref{eq:cota-defecto-potencias-cadena-bochner-cap13} give
\[
\begin{split}
\|\nabla_\ell \mathbf{z}\|_{\mathscr D_{\ell+1,j}^p}
&\leq
C\|B_{\ell,p}^{\frac{j}{2}}\mathbf{z}\|_{\mathscr D_{\ell,1}^p}
+
C\|\mathbf{z}\|_{\mathscr D_{\ell,j-1}^p}\\
&\leq
C'\|\mathbf{z}\|_{\mathscr D_{\ell,j+1}^p}.
\end{split}
\]
We have therefore proved the continuous mapping
\begin{equation}
\label{eq:gradiente-dominios-cadena-bochner-cap13}
\nabla_\ell\colon
\mathscr D_{\ell,j+1}^p
\longrightarrow
\mathscr D_{\ell+1,j}^p.
\end{equation}

We now prove \eqref{eq:estimaciones-meyer-cadena-bochner-cap13} by induction on $k$. The case $k=0$ is immediate, and the case $k=1$ is contained in \eqref{eq:meyer-primer-orden-cadena-bochner-cap13} and \eqref{eq:cotas-meyer-primer-orden-cadena-bochner-cap13}. Fix $k\in\{2,\dots,N\}$, suppose the assertion has been proved at all orders less than $k$, and fix
\[
0\leq\ell\leq N-k.
\]
For $\rho\in\{0,\dots,k-2\}$ and $a\in\{0,\dots,\rho\}$, the Leibniz rule and \eqref{eq:cotas-coeficientes-cadena-bochner-cap13} give
\begin{equation}
\label{eq:multiplicacion-W-cadena-bochner-cap13}
\|\nabla^a(\mathbf{W}_\ell \mathbf{z})\|_{L^p(M,\mathbf{E}_{\ell+a+1})}
\leq
C_a\displaystyle\sum_{b=0}^a
\|\nabla^b \mathbf{z}\|_{L^p(M,\mathbf{E}_{\ell+b})}.
\end{equation}
Summing over $a$ and applying the induction hypothesis at order $\rho$ yields a continuous operator
\begin{equation}
\label{eq:W-en-dominios-cadena-bochner-cap13}
\mathbf{W}_\ell\colon
\mathscr D_{\ell,\rho}^p
\longrightarrow
\mathscr D_{\ell+1,\rho}^p,
\qquad
0\leq\rho\leq k-2.
\end{equation}
The indices are admissible because $\ell+1\leq N-(k-2)$. Consequently, we may use \eqref{eq:gradiente-dominios-cadena-bochner-cap13} and \eqref{eq:cota-defecto-potencias-cadena-bochner-cap13} with $j=k-1$.

Let $\mathbf{u}\in\mathscr D_{\ell,k}^p$. By \eqref{eq:gradiente-dominios-cadena-bochner-cap13},
\[
\nabla_\ell \mathbf{u}\in\mathscr D_{\ell+1,k-1}^p.
\]
The inductive estimate at order $k-1$, applied to this section, and the identity defining $\Theta_{\ell,k-1}$ give
\begin{equation*}
\begin{split}
\|\nabla^k \mathbf{u}\|_{L^p(M,\mathbf{E}_{\ell+k})}
&\leq
C\|B_{\ell+1,p}^{\frac{k-1}{2}}\nabla_\ell \mathbf{u}\|_{L^p(M,\mathbf{E}_{\ell+1})}\\
&\leq
C\|\nabla_\ell B_{\ell,p}^{\frac{k-1}{2}}\mathbf{u}\|_{L^p(M,\mathbf{E}_{\ell+1})}
+C\|\Theta_{\ell,k-1}\mathbf{u}\|_{L^p(M,\mathbf{E}_{\ell+1})}\\
&\leq
C\|B_{\ell,p}^{\frac{k}{2}}\mathbf{u}\|_{L^p(M,\mathbf{E}_\ell)}
+C\|B_{\ell,p}^{\frac{k-2}{2}}\mathbf{u}\|_{L^p(M,\mathbf{E}_\ell)}\\
&\leq
C'\|B_{\ell,p}^{\frac{k}{2}}\mathbf{u}\|_{L^p(M,\mathbf{E}_\ell)}.
\end{split}
\end{equation*}
In the third line, we use \eqref{eq:cotas-meyer-primer-orden-cadena-bochner-cap13} for the section $B_{\ell,p}^{\frac{k-1}{2}}\mathbf{u}$ and \eqref{eq:cota-defecto-potencias-cadena-bochner-cap13}; in the last, we use \eqref{eq:monotonia-dominios-bessel-cadena-cap13}. For lower-order derivatives, the induction hypothesis and the same monotonicity give
\[
\displaystyle\sum_{j=0}^{k-1}
\|\nabla^j \mathbf{u}\|_{L^p(M,\mathbf{E}_{\ell+j})}
\leq
C\|B_{\ell,p}^{\frac{k-1}{2}}\mathbf{u}\|_{L^p(M,\mathbf{E}_\ell)}
\leq
C\|B_{\ell,p}^{\frac{k}{2}}\mathbf{u}\|_{L^p(M,\mathbf{E}_\ell)}.
\]
Therefore,
\begin{equation}
\label{eq:meyer-derivadas-por-potencia-cap13}
\displaystyle\sum_{j=0}^k
\|\nabla^j \mathbf{u}\|_{L^p(M,\mathbf{E}_{\ell+j})}
\leq
C\|B_{\ell,p}^{\frac{k}{2}}\mathbf{u}\|_{L^p(M,\mathbf{E}_\ell)}.
\end{equation}
In particular, this estimate proves the continuous inclusion
$\mathscr D_{\ell,k}^p\subseteq W^{k,p}(M,\mathbf{E}_\ell)$

We prove the reverse estimate on smooth compactly supported sections. Proposition~\ref{prop:generador-Lp-extiende-laplaciano} shows that $\Gamma_c(\mathbf{E}_\ell)\subseteq\mathcal D(L_{\ell,p})$. Since the differential operator $L_{\ell,0}+\mathbf{V}_\ell$ preserves $\Gamma_c(\mathbf{E}_\ell)$, $B_{\ell,p}$ preserves it as well. For $r=1$, we already have $\Gamma_c(\mathbf E_\ell)\subseteq\mathcal D(B_{\ell,p})$. If $\mathbf u\in\Gamma_c(\mathbf E_\ell)$ belongs to $\mathcal D(B_{\ell,p}^r)$, then $B_{\ell,p}^r\mathbf u\in\Gamma_c(\mathbf E_\ell)
\subseteq\mathcal D(B_{\ell,p})$; by definition of the domain of a composition, $\mathbf u\in\mathcal D(B_{\ell,p}^{r+1})$. Induction gives
\[
\Gamma_c(\mathbf{E}_\ell)
\subseteq
\bigcap_{r\in\mathbb N}\mathcal D(B_{\ell,p}^r).
\]
Monotonicity of the scale implies that every section in $\Gamma_c(\mathbf{E}_\ell)$ also belongs to the domains of the fractional powers appearing below. For $\mathbf{u}\in\Gamma_c(\mathbf{E}_\ell)$, apply the first-order estimate to $B_{\ell,p}^{\frac{k-1}{2}}\mathbf{u}$ and then add and subtract the defect $\Theta_{\ell,k-1}\mathbf{u}$:
\begin{equation*}
\begin{split}
\|B_{\ell,p}^{\frac{k}{2}}\mathbf{u}\|_{L^p(M,\mathbf{E}_\ell)}
&\leq
C\|B_{\ell,p}^{\frac{k-1}{2}}\mathbf{u}\|_{L^p(M,\mathbf{E}_\ell)}
+C\|\nabla_\ell B_{\ell,p}^{\frac{k-1}{2}}\mathbf{u}\|_{L^p(M,\mathbf{E}_{\ell+1})}\\
&\leq
C\|B_{\ell,p}^{\frac{k-1}{2}}\mathbf{u}\|_{L^p(M,\mathbf{E}_\ell)}
+C\|B_{\ell+1,p}^{\frac{k-1}{2}}\nabla_\ell \mathbf{u}\|_{L^p(M,\mathbf{E}_{\ell+1})}\\
&\hspace{12mm}
+C\|\Theta_{\ell,k-1}\mathbf{u}\|_{L^p(M,\mathbf{E}_{\ell+1})}\\
&\leq
C'\displaystyle\sum_{j=0}^k
\|\nabla^j \mathbf{u}\|_{L^p(M,\mathbf{E}_{\ell+j})}.
\end{split}
\end{equation*}
In the last inequality, apply the induction hypothesis at order $k-1$ to $\mathbf{u}$ and $\nabla_\ell \mathbf{u}$, and at order $k-2$ to the term controlling $\Theta_{\ell,k-1}\mathbf{u}$. Thus,
\begin{equation}
\label{eq:meyer-potencia-por-derivadas-nucleo-cap13}
\|B_{\ell,p}^{\frac{k}{2}}\mathbf{u}\|_{L^p(M,\mathbf{E}_\ell)}
\leq
C\displaystyle\sum_{j=0}^k
\|\nabla^j \mathbf{u}\|_{L^p(M,\mathbf{E}_{\ell+j})},
\qquad
\mathbf{u}\in\Gamma_c(\mathbf{E}_\ell).
\end{equation}
Let us verify the required density. Let $\mathcal T^{\operatorname{geo}}=(U_i,\phi_i,h_i)_{i\in I}$ be a geodesic trivialization, so that $\displaystyle\sum_{i\in I}h_i=1$. For finite $F\subseteq I$, set
\[
\mathbf{u}_F:=\displaystyle\sum_{i\in F}h_i \mathbf{u}.
\]
Each $\mathbf{u}_F$ has compact support. The same Leibniz calculation and second overlap count used in Proposition~\ref{prop:covariantes-localizacion-orden-entero} give
\begin{equation}
\label{eq:cola-densidad-Wkp-geometria-acotada-cap13}
\|\mathbf{u}-\mathbf{u}_F\|_{W^{k,p}(M,\mathbf{E}_\ell)}^p
\leq
C\displaystyle\sum_{i\in I\setminus F}
\|h_i \mathbf{u}\|_{W^{k,p}(U_i,\mathbf{E}_\ell)}^p.
\end{equation}
The series on the right converges by the direct localization inequality, so its tails tend to zero. Finally, the Meyers--Serrin Theorem~\ref{meyers-serrin-haz} approximates each $\mathbf{u}_F$ by smooth sections. If $\chi_F\in C_c^\infty(M)$ equals one on a neighborhood of $\operatorname{supp}(\mathbf{u}_F)$, multiplying those approximations by $\chi_F$ produces a sequence in $\Gamma_c(\mathbf{E}_\ell)$ converging to $\mathbf{u}_F$ in $W^{k,p}(M,\mathbf E_\ell)$. Therefore, $\Gamma_c(\mathbf{E}_\ell)$ is dense in $W^{k,p}(M,\mathbf{E}_\ell)$.

Now let $\mathbf{u}\in W^{k,p}(M,\mathbf{E}_\ell)$, and choose $\mathbf{u}_n\in\Gamma_c(\mathbf{E}_\ell)$ with $\mathbf{u}_n\to \mathbf{u}$ in $W^{k,p}(M,\mathbf E_\ell)$. Estimate \eqref{eq:meyer-potencia-por-derivadas-nucleo-cap13} shows that $(B_{\ell,p}^{\frac{k}{2}}\mathbf{u}_n)_n$ is Cauchy in $L^p(M,\mathbf{E}_\ell)$. Moreover, $\mathbf{u}_n\to \mathbf{u}$ in $L^p(M,\mathbf{E}_\ell)$. Positive powers in the Bessel scale are closed operators; consequently,
\[
\mathbf{u}\in\mathscr D_{\ell,k}^p
\quad\text{and}\quad
B_{\ell,p}^{\frac{k}{2}}\mathbf{u}
=
\lim_{n\to\infty}B_{\ell,p}^{\frac{k}{2}}\mathbf{u}_n.
\]
Passing to the limit gives
\begin{equation}
\label{eq:meyer-potencia-por-derivadas-cap13}
\|B_{\ell,p}^{\frac{k}{2}}\mathbf{u}\|_{L^p(M,\mathbf{E}_\ell)}
\leq
C\displaystyle\sum_{j=0}^k
\|\nabla^j \mathbf{u}\|_{L^p(M,\mathbf{E}_{\ell+j})}.
\end{equation}
Together with \eqref{eq:meyer-derivadas-por-potencia-cap13}, this proves equality of spaces and both estimates at order $k$. The induction is complete.
\end{proof}

\begin{theorem}[Bessel potentials and covariant derivatives]
\label{teo:yoshida-bessel-covariante-geometria-acotada}
Let $(M,\mathbf{g})$ be a Riemannian manifold without boundary and with bounded geometry, and let $(\mathbf{E},\mathbf{h}_{\mathbf{E}},\nabla^{\mathbf{E}})$ be a bundle of bounded geometry. If $m\in\mathbb N_0$ and $1<p<\infty$, then
\[
H_{L_{\mathbf{E}}}^{m,p}(M,\mathbf{E})=W^{m,p}(M,\mathbf{E}).
\]
There exist constants $c_{m,p},C_{m,p}>0$, independent of $\mathbf{u}$, such that
\begin{equation}
\label{eq:yoshida-bessel-covariante}
c_{m,p}\|\mathbf{u}\|_{W^{m,p}(M,\mathbf{E})}
\leq
\|\mathbf{u}\|_{H_{L_{\mathbf{E}}}^{m,p}(M,\mathbf{E})}
\leq
C_{m,p}\|\mathbf{u}\|_{W^{m,p}(M,\mathbf{E})}.
\end{equation}
\end{theorem}

\begin{proof}
The case $m=0$ is immediate. Suppose $m\geq1$ and apply Proposition~\ref{prop:estimaciones-meyer-cadena-bochner-geometria-acotada} with $N=m$, $\ell=0$, and $k=m$. Since $\mathbf{V}_0=0$, there exist $a_{m,p},b_{m,p}>0$ such that
\begin{equation}
\label{eq:meyer-desplazamiento-mu-cap13}
a_{m,p}
\displaystyle\sum_{j=0}^m
\|\nabla^j \mathbf{u}\|_{L^p(M,T^{(0,j)}(TM)\otimes \mathbf{E})}
\leq
\|(\mu_m I+L_{\mathbf{E},p})^{\frac{m}{2}}\mathbf{u}\|_{L^p(M,\mathbf{E})}
\leq
b_{m,p}
\displaystyle\sum_{j=0}^m
\|\nabla^j \mathbf{u}\|_{L^p(M,T^{(0,j)}(TM)\otimes \mathbf{E})}.
\end{equation}
Lemma~\ref{lem:independencia-desplazamiento-bessel-Lp}, applied to the covariant semigroup, shows that
\[
\operatorname{Ran}\bigl((\mu_m I+L_{\mathbf{E},p})^{-\frac{m}{2}}\bigr)
=
\operatorname{Ran}\bigl((I+L_{\mathbf{E},p})^{-\frac{m}{2}}\bigr)
=
H_{L_{\mathbf{E}}}^{m,p}(M,\mathbf{E}),
\]
and, for $\mathbf{u}$ in this common range,
\[
 \|Q_{\mu_m,1}\|_{\mathcal L(L^p(M,\mathbf{E}))}^{-1}
 \|\mathbf{u}\|_{\operatorname{Ran}((I+L_{\mathbf{E},p})^{-\frac{m}{2}})}
 \leq
 \|\mathbf{u}\|_{\operatorname{Ran}((\mu_mI+L_{\mathbf{E},p})^{-\frac{m}{2}})}
 \leq
 \|Q_{1,\mu_m}\|_{\mathcal L(L^p(M,\mathbf{E}))}
 \|\mathbf{u}\|_{\operatorname{Ran}((I+L_{\mathbf{E},p})^{-\frac{m}{2}})},
\]
where the operators $Q_{a,b}$ are those constructed in the proof of that lemma. On the other hand, the finite-sum norms $\ell^p(\{0,\ldots,k\})$ and $\ell^1(\{0,\ldots,k\})$ satisfy
\begin{equation}
\label{eq:suma-l1-lp-derivadas-cap13}
\begin{split}
\left(
\displaystyle\sum_{j=0}^m
\|\nabla^j \mathbf{u}\|_{L^p(M,T^{(0,j)}(TM)\otimes \mathbf{E})}^p
\right)^{\frac1p}
&\leq
\displaystyle\sum_{j=0}^m
\|\nabla^j \mathbf{u}\|_{L^p(M,T^{(0,j)}(TM)\otimes \mathbf{E})}\\
&\leq
(m+1)^{1-\frac1p}
\left(
\displaystyle\sum_{j=0}^m
\|\nabla^j \mathbf{u}\|_{L^p(M,T^{(0,j)}(TM)\otimes \mathbf{E})}^p
\right)^{\frac1p}.
\end{split}
\end{equation}
The parenthesized expression is the norm of $W^{m,p}(M,\mathbf{E})$. The preceding comparisons prove equality of spaces and \eqref{eq:yoshida-bessel-covariante}.

In the Hermitian case, the preceding estimates are interpreted on the underlying real bundle when the Bellman function is applied. The connection, potentials $\mathbf V_\ell$, operators $B_{\ell,p}$, and covariant derivatives commute with the complex structure, and the norm of the underlying real bundle agrees with the Hermitian norm. Theorem~\ref{teo:potencias-imaginarias-bochner-schrodinger-autocontenido} also proves that the extensions obtained are complex-linear. Thus, the preceding comparisons retain the same constants in the complex case.
\end{proof}

To interpolate the scale defined through heat, we need to control powers of imaginary order. Self-adjointness resolves this on $L^2(M,\mathbf E)$. For other exponents, we can represent a section by scalar functions on the unit spheres of its fibers. This representation preserves the $L^p(M,\mathbf E)$ norm up to a fixed factor and converts the covariant Laplacian into a shifted scalar Laplacian.

\begin{lemma}[Imaginary powers of the covariant Laplacian]
\label{lem:potencias-imaginarias-bochner-haz-completo}
Let $(M,\mathbf g)$ be a complete Riemannian manifold without boundary, and let $\mathbf E\to M$ be a finite-rank vector bundle equipped with a bundle metric $\mathbf h_{\mathbf E}$, Hermitian in the complex case, and a compatible connection $\nabla^{\mathbf E}$. For $1<p<\infty$, let $L_{\mathbf E,p}$ be the positive generator of covariant heat and set $B_p:=I+L_{\mathbf E,p}$. The spectral powers $B_2^{it}$, with $t\in\mathbb R$, extend consistently to bounded operators on $L^p(M,\mathbf E)$. There exist $C_{p,\mathbf E},\omega_{p,\mathbf E}>0$, depending only on $p$, the rank, and the scalar field, such that
\[
 \|B_p^{it}\|_{\mathcal L(L^p(M,\mathbf E))}
 \leq C_{p,\mathbf E}e^{\omega_{p,\mathbf E}|t|},
 \qquad t\in\mathbb R.
\]
In the real case, these powers are understood on the complexification.
\end{lemma}

\begin{proof}
If the bundle is real, complexify it with the induced Hermitian metric and connection. The original norm is preserved on real sections, and the complex estimates restrict to them. It therefore suffices to treat a complex bundle with a Hermitian bundle metric. Write $r=\operatorname{rank}_{\mathbb C}\mathbf E$ and $d=2r$. Consider the manifold
\[
 \widehat M=\mathbb S(\mathbf E)
 :=\{(x,\mathbf v)\in\mathbf E\mid |\mathbf v|_{\mathbf h}=1\},
 \qquad \pi(x,\mathbf v)=x.
\]
The connection splits its tangent space into horizontal and vertical parts. Equip the horizontal part with the metric of $M$, each fiber with the unit-sphere metric, and declare the two parts orthogonal. Denote the resulting metric by $\widehat{\mathbf g}$.

This metric is complete. Indeed, $\pi$ does not increase lengths. A Cauchy sequence in $\widehat M$ projects to a Cauchy sequence in $M$, converging to some $x$. Apart from finitely many terms, it lies over a compact set contained in a bundle-trivializing neighborhood of $x$. The inverse image of that compact set is compact because the fiber is a sphere. Thus, the original sequence has a convergent subsequence and, being Cauchy, converges. The same argument applies componentwise if the sphere bundle is disconnected.

Let $d\sigma_x$ be the invariant probability measure on the sphere of $\mathbf E_x$. The volume measure of $\widehat{\mathbf g}$, divided by the area of $\mathbb S^{d-1}$, is $d\lambda_{\mathbf g}(x)\,d\sigma_x(\mathbf v)$; the preceding orthogonal decomposition verifies this in any local orthonormal frame. Use this normalized measure on $L^p(\widehat M)$. Normalization does not change operator norms, since it multiplies the norms of arguments and images by the same factor. Define
\[
 (\mathcal J\mathbf u)(x,\mathbf v)
 :=\langle\mathbf u(x),\mathbf v\rangle_{\mathbf h}.
\]
Use the Hermitian product linear in its first variable. Invariance of $\sigma_x$ under fiber isometries gives
\begin{equation}
\label{eq:norma-elevacion-esferas-haz}
 \|\mathcal J\mathbf u\|_{L^p(\widehat M)}
 =c_{r,p}\|\mathbf u\|_{L^p(M,\mathbf E)},
 \qquad
 c_{r,p}^{p}
 :=\int_{\mathbb S^{d-1}}
       |\langle\mathbf e_1,\mathbf v\rangle|^p\,d\sigma(\mathbf v)>0.
\end{equation}
Consequently, $\mathcal J$ has closed image. For $p=2$, we have $c_{r,2}=r^{-1/2}$.

We verify the relation between the Laplacians. If $\mathbf X$ is a field on $M$ and $\mathbf X^H$ is its horizontal lift, the definition of parallel transport and metric compatibility imply
\[
 \mathbf X^H(\mathcal J\mathbf u)
 =\mathcal J(\nabla_{\mathbf X}\mathbf u).
\]
The fibers are totally geodesic for $\widehat{\mathbf g}$. To see this, transport two vertical fields $\mathbf V,\mathbf W$ by the horizontal flow of $\mathbf X^H$. This transport is an isometry of the spheres, so
\[
 \mathbf X^H\langle\mathbf V,\mathbf W\rangle
 =\langle[\mathbf X^H,\mathbf V],\mathbf W\rangle
  +\langle\mathbf V,[\mathbf X^H,\mathbf W]\rangle.
\]
Substituting this identity into the Koszul formula gives $\langle\widehat\nabla_{\mathbf V}\mathbf W,\mathbf X^H\rangle=0$. For horizontal fields, every Koszul term paired with $\mathbf Z^H$ is the lift of the corresponding term on $M$. Hence,
\[
 \langle\widehat\nabla_{\mathbf X^H}\mathbf Y^H,\mathbf Z^H\rangle
 =\mathbf g(\nabla^M_{\mathbf X}\mathbf Y,\mathbf Z)\circ\pi.
\]
The vertical component of the repeated derivative is computed with the same formula:
\[
 2\langle\widehat\nabla_{\mathbf X^H}\mathbf X^H,\mathbf V\rangle
 =-\mathbf V\bigl(\mathbf g(\mathbf X,\mathbf X)\circ\pi\bigr)
   -2\langle[\mathbf X^H,\mathbf V],\mathbf X^H\rangle=0.
\]
The first term is zero because $\mathbf V$ is vertical; the second is zero because the bracket of a projectable field with a vertical field is vertical. Thus, taking the horizontal trace of the second derivative of $\mathcal J\mathbf u$ gives exactly $\mathcal J(\operatorname{tr}_{\mathbf g}\nabla^2\mathbf u)$.

For the vertical trace, fix $\mathbf v$ in a fiber and a vector $\mathbf w$ tangent to its sphere. The spherical geodesic starting at $\mathbf v$ with velocity $\mathbf w$ satisfies $\boldsymbol\gamma''(0)=-|\mathbf w|^2\mathbf v$. The function $\langle\mathbf u(x),\boldsymbol\gamma(\tau)\rangle$ therefore has second derivative at zero equal to $-|\mathbf w|^2\langle\mathbf u(x),\mathbf v\rangle$. Take an orthonormal basis of the $d-1$ tangent vectors of the fiber and sum: the positive-sign vertical Laplacian is $(d-1)$ times $\mathcal J\mathbf u$. Adding the horizontal and vertical traces gives, for $\mathbf u\in\Gamma_c(\mathbf E)$,
\begin{equation}
\label{eq:laplaciano-esferas-lineales-haz}
 L_{\widehat M}\mathcal J\mathbf u
 =\mathcal J\bigl(L_{\mathbf E,0}\mathbf u+(d-1)\mathbf u\bigr).
\end{equation}

The differential identity also determines the closed realizations. Set $A:=I+L_{\widehat M}$ and $D_2:=dI+L_{\mathbf E,2}$. Completeness and the essential self-adjointness theorem of Chapter~\ref{cap:sobolev-fraccionario-haces-nucleo-calor} allow each $\mathbf u\in\mathcal D(D_2)$ to be approximated by sections in $\Gamma_c(\mathbf E)$ in the graph norm. Their images under $\mathcal J$ have compact support because $\pi$ has compact fibers. Closedness of $A$ and \eqref{eq:laplaciano-esferas-lineales-haz} give
\[
 \mathcal J\mathcal D(D_2)\subseteq\mathcal D(A),
 \qquad A\mathcal J\mathbf u=\mathcal JD_2\mathbf u.
\]
If $\lambda>0$ and $\mathbf f\in L^2(M,\mathbf E)$, apply this equality to $\mathbf u=(\lambda+D_2)^{-1}\mathbf f$. Uniqueness of the solution of the resolvent equation implies
\[
 (\lambda+A)^{-1}\mathcal J\mathbf f
 =\mathcal J(\lambda+D_2)^{-1}\mathbf f.
\]
The resolvents are self-adjoint; the closed image of $\mathcal J$ and its orthogonal complement are invariant. Spectral calculus on that subspace then gives
$A^{it}\mathcal J=\mathcal JD_2^{it}$
for $t\in\mathbb R$.

The scalar case of Theorem~\ref{teo:potencias-imaginarias-bochner-schrodinger-autocontenido}, applied to the operator $A=I+L_{\widehat M,p}$ on the complete manifold $(\widehat M,\widehat{\mathbf g})$, gives
\[
 \|A^{it}\|_{\mathcal L(L^p(\widehat M))}
 \leq M_p e^{\omega_p|t|},\qquad t\in\mathbb R.
\]
For $\mathbf u\in L^2(M,\mathbf E)\cap L^p(M,\mathbf E)$, use spectral intertwining and \eqref{eq:norma-elevacion-esferas-haz}; the factor $c_{r,p}$ appears on both sides and cancels. This yields
\[
 \|D_2^{it}\mathbf u\|_{L^p(M,\mathbf E)}
 \leq M_p e^{\omega_p|t|}\|\mathbf u\|_{L^p(M,\mathbf E)}.
\]
Density of the intersection defines the extension $D_p^{it}$ and ensures its consistency, where $D_p=dI+L_{\mathbf E,p}$.

Finally, pass from the shift $d$ to the shift $1$. The integral representation of the resolvent and heat contractivity give $\|D_p^{-1}\|\leq d^{-1}$. Set
\[
 K_p:=(d-1)D_p^{-1},\qquad
 G_p:=-\displaystyle\sum_{j=1}^{\infty}\frac{K_p^j}{j}.
\]
This series converges in operator norm, since $\|K_p\|\leq(d-1)/d<1$, and
\[
 \|G_p\|\leq
 \displaystyle\sum_{j=1}^{\infty}\frac{((d-1)/d)^j}{j}=\log d.
\]
On $L^2(M,\mathbf E)$, $G_2=\log(I-K_2)$ and $B_2=D_2(I-K_2)$, so $B_2^{it}=D_2^{it}e^{itG_2}$. The series are consistent on the spaces $L^p(M,\mathbf E)$, and hence
\[
 \|B_p^{it}\|
 \leq M_pe^{\omega_p|t|}e^{|t|\|G_p\|}
 \leq M_pe^{(\omega_p+\log d)|t|}.
\]
We also obtain strong continuity in the parameter. If $p\neq2$, choose $q\in(1,\infty)$ so that $p$ lies strictly between $2$ and $q$, and $\theta\in(0,1)$ with $1/p=\theta/2+(1-\theta)/q$. For $\mathbf u\in\Gamma_c(\mathbf E)$, Hölder gives
\[
 \|(B^{it}-B^{it_0})\mathbf u\|_{L^{p}(M,\mathbf E)}
 \leq\|(B^{it}-B^{it_0})\mathbf u\|_{L^{2}(M,\mathbf E)}^\theta
      \|(B^{it}-B^{it_0})\mathbf u\|_{L^{q}(M,\mathbf E)}^{1-\theta}.
\]
The first factor tends to zero by spectral calculus; the second is bounded when $t$ remains in a compact interval, by the estimate already proved with exponent $q$. Density of $\Gamma_c(\mathbf E)$ in $L^p(M,\mathbf E)$ and the local uniform operator bound extend convergence to every section in $L^p(M,\mathbf E)$. For $p=2$, spectral continuity applies directly.
\end{proof}

\begin{theorem}[Interpolation of Bessel potentials in bundles]
\label{teo:interpolacion-compleja-bessel-laplaciano-completo}
Let $(M,\mathbf{g})$ be a complete Riemannian manifold without boundary, and let $\mathbf E\to M$ be a finite-rank vector bundle equipped with a bundle metric $\mathbf h_{\mathbf E}$, Hermitian in the complex case, and a compatible connection $\nabla^{\mathbf E}$. Let $L_{\mathbf E,p}$ be the positive generator of covariant heat and let $1<p<\infty$. For real bundles, use complexification and restrict the interpolated spaces to real sections. If $0\leq s_0<s_1$, $0<\theta<1$, and $s=(1-\theta)s_0+\theta s_1$, then
\begin{equation}
\label{eq:interpolacion-compleja-bessel-intrinseco}
[H_{L_{\mathbf E}}^{s_0,p}(M,\mathbf E),H_{L_{\mathbf E}}^{s_1,p}(M,\mathbf E)]_\theta
=
H_{L_{\mathbf E}}^{s,p}(M,\mathbf E)
\end{equation}
and there exist constants $c_{s_0,s_1,\theta,p},C_{s_0,s_1,\theta,p}>0$ such that
\[
 c_{s_0,s_1,\theta,p}\|\mathbf{u}\|_{H_{L_{\mathbf E}}^{s,p}(M,\mathbf E)}
 \leq
 \|\mathbf{u}\|_{[H_{L_{\mathbf E}}^{s_0,p}(M,\mathbf E),H_{L_{\mathbf E}}^{s_1,p}(M,\mathbf E)]_\theta}
 \leq
C_{s_0,s_1,\theta,p}\|\mathbf{u}\|_{H_{L_{\mathbf E}}^{s,p}(M,\mathbf E)}.
\]
The constants may also depend on the rank of the bundle, which is fixed.
\end{theorem}

\begin{proof}
Set
\[
B:=I+L_{\mathbf E,p},
\qquad
a_j:=\frac{s_j}{2}\quad(j\in\{0,1\}),
\qquad
a:=(1-\theta)a_0+\theta a_1,
\qquad
\delta:=a_1-a_0.
\]
For $r\geq0$, the potential $B^{-r}$ is injective and its range is the domain of $B^r$; moreover,
\[
\|\mathbf{u}\|_{H_{L_{\mathbf E}}^{2r,p}(M,\mathbf E)}
=
\|B^r\mathbf{u}\|_{L^p(M,\mathbf E)}.
\]
It therefore suffices to interpolate the domains $\mathcal D(B^{a_j})$ equipped with these norms.

Lemma~\ref{lem:potencias-imaginarias-bochner-haz-completo} proves that imaginary powers of $B$ are bounded on $L^p(M,\mathbf E)$ and that there exist $M_{p,\mathbf E},\omega_{p,\mathbf E}>0$ such that
\begin{equation}
\label{eq:potencias-imaginarias-strichartz-cap13}
\|B^{it}\|_{\mathcal L(L^p(M,\mathbf E))}
\leq
M_{p,\mathbf E} e^{\omega_{p,\mathbf E}|t|},
\qquad t\in\mathbb R.
\end{equation}
This is the analytic estimate allowing the interpolation argument of Chapter~\ref{cap:espacios-interpolacion} to be transferred to the manifold. Use its strip $\mathcal T=\{z\in\mathbb C\mid0<\operatorname{Re}z<1\}$.

Fix an integer $N>a_1$ and first consider $\mathbf{u}\in\mathcal D(B^N)$. For $\varepsilon>0$ and $z\in\overline{\mathcal T}$, set
$a(z):=(1-z)a_0+za_1$
and
\[
F_\varepsilon(z)
:=
e^{\varepsilon(z-\theta)^2}B^{a-a(z)}\mathbf{u}.
\]
Holomorphy is verified by writing $B^{a-a(z)}\mathbf u=B^{a-a(z)-N}B^N\mathbf u$: the exponent of the first factor has negative real part throughout the closed strip, and its Laplace integral is holomorphic in the interior. For $z=\sigma+it$, with $\sigma\in[0,1]$, it also factors into the negative real power and $B^{-it\delta}$. The former has norm at most one, and the latter satisfies \eqref{eq:potencias-imaginarias-strichartz-cap13}. The Gaussian factor makes the family bounded on the strip and its boundary values tend to zero. We have $F_\varepsilon(\theta)=\mathbf u$. If $j\in\{0,1\}$ and $t\in\mathbb R$, then
\begin{equation}
\label{eq:frontera-familia-potencias-cap13}
B^{a_j}F_\varepsilon(j+it)
=
e^{\varepsilon(j+it-\theta)^2}
B^{-it\delta}B^a\mathbf{u}.
\end{equation}
By \eqref{eq:potencias-imaginarias-strichartz-cap13}, the norm of the right-hand side is bounded by
\[
M_{p,\mathbf E}
\exp\bigl(
\varepsilon(j-\theta)^2-
\varepsilon t^2+
\omega_{p,\mathbf E}\delta|t|
\bigr)
\|B^a\mathbf{u}\|_{L^p(M,\mathbf E)}.
\]
Continuity in each endpoint norm follows from the same factorization and strong continuity of the imaginary powers. This expression is bounded and tends to zero as $|t|\to\infty$. Thus, $F_\varepsilon$ belongs to the Calderón space of the pair $(\mathcal D(B^{a_0}),\mathcal D(B^{a_1}))$, and
\begin{equation}
\label{eq:inclusion-dominio-interpolado-cap13}
\|\mathbf{u}\|_{[\mathcal D(B^{a_0}),\mathcal D(B^{a_1})]_\theta}
\leq
C\|B^a\mathbf{u}\|_{L^p(M,\mathbf E)}.
\end{equation}
Now let $\mathbf{u}\in\mathcal D(B^a)$ and set $\mathbf{u}_\eta=(I+\eta B)^{-N}\mathbf{u}$. Then $\mathbf{u}_\eta\in\mathcal D(B^N)$ and, since the resolvent commutes with $B^a$, we have $B^a\mathbf{u}_\eta\to B^a\mathbf{u}$ in $L^p(M,\mathbf E)$. Estimate \eqref{eq:inclusion-dominio-interpolado-cap13}, applied to $\mathbf{u}_\eta-\mathbf{u}_{\eta\prime}$, shows that $(\mathbf{u}_\eta)$ is Cauchy in the interpolated space. Its limit agrees with $\mathbf{u}$ in $L^p(M,\mathbf E)$ because the inclusions of both endpoint spaces in $L^p(M,\mathbf E)$ are continuous. This proves the inclusion and estimate for every $\mathbf{u}\in\mathcal D(B^a)$.

We prove the reverse inclusion. Let $F\in\mathscr F(\mathcal D(B^{a_0}),\mathcal D(B^{a_1}))$ and write $\mathbf{u}=F(\theta)$. Regularize the values of $F$ using operators that map them into the domains of the powers we will use. Fix an integer $N>a_1$, take $\eta>0$, and set
\[
R_\eta:=(I+\eta B)^{-N}.
\]
The operator $R_\eta$ commutes with powers of $B$ and maps $L^p(M,\mathbf E)$ into the domain of $B^N$. Thus, the function
\[
G_{\varepsilon,\eta}(z)
:=
e^{\varepsilon(z-\theta)^2}
B^{a(z)}R_\eta F(z)
\]
is holomorphic with values in $L^p(M,\mathbf E)$. On the boundary,
\[
G_{\varepsilon,\eta}(j+it)
=
e^{\varepsilon(j+it-\theta)^2}
B^{it\delta}R_\eta B^{a_j}F(j+it).
\]
Uniform contractivity of $R_\eta$, the bound \eqref{eq:potencias-imaginarias-strichartz-cap13}, and the Gaussian factor give a boundary bound independent of $\eta$. The three-lines lemma, applied to the pair $(L^p(M,\mathbf E),L^p(M,\mathbf E))$, gives
\begin{equation}
\label{eq:cota-potencia-regularizada-cap13}
\|B^aR_\eta \mathbf{u}\|_{L^p(M,\mathbf E)}
\leq
C\|F\|_{\mathscr F}.
\end{equation}
As $\eta\to0^+$, $R_\eta \mathbf{u}\to \mathbf{u}$ in $L^p(M,\mathbf E)$. Since $1<p<\infty$, from \eqref{eq:cota-potencia-regularizada-cap13} we can extract a subsequence along which $B^aR_\eta \mathbf{u}$ converges weakly in $L^p(M,\mathbf E)$. The graph of the closed operator $B^a$ is a closed subspace and hence weakly closed. It follows that $\mathbf{u}\in\mathcal D(B^a)$ and
\[
\|B^a\mathbf{u}\|_{L^p(M,\mathbf E)}
\leq
C\|F\|_{\mathscr F}.
\]
Taking the infimum over functions $F$ with $F(\theta)=\mathbf{u}$ proves the reverse inclusion in \eqref{eq:inclusion-dominio-interpolado-cap13}. The constant in that estimate and the constant in \eqref{eq:inclusion-dominio-interpolado-cap13} are $c_{s_0,s_1,\theta,p}^{-1}$ and $C_{s_0,s_1,\theta,p}$, respectively.

Lemma~\ref{lem:independencia-desplazamiento-bessel-Lp} allows any positive shift to be replaced by the shift $1$ chosen in the definition of $B$.
\end{proof}

\begin{theorem}[Uniform localization and Bessel potentials]
\label{teo:triebel-bessel-localizacion-geometria-acotada}
Let $(M,\mathbf{g})$ be a Riemannian manifold without boundary and with bounded geometry, let $s\in\mathbb R$, and let $1<p<\infty$. If $L=\Delta_B=-\Delta$ is the positive realization of the scalar Laplacian, then, as spaces of distributions,
\[
H^{s,p}(M)=F^s_{p,2}(M)=H_L^{s,p}(M).
\]
More precisely, there exist constants $c_{s,p},C_{s,p}>0$ such that
\begin{equation}
\label{eq:equivalencia-local-intrinseca-bessel}
c_{s,p}\|u\|_{H^{s,p}(M)}
\leq
\|u\|_{H_L^{s,p}(M)}
\leq
C_{s,p}\|u\|_{H^{s,p}(M)},
\end{equation}
and there exist $d_{s,p},D_{s,p}>0$ such that
\begin{equation}
\label{eq:equivalencia-local-H-F}
d_{s,p}\|u\|_{F^s_{p,2}(M)}
\leq
\|u\|_{H^{s,p}(M)}
\leq
D_{s,p}\|u\|_{F^s_{p,2}(M)}.
\end{equation}
The constants do not depend on $u$.
\end{theorem}

The interpolation scheme for this identification is that of \cite[Section~7.4.5]{Triebel2}.

\begin{proof}
Begin with the identification of the two localized scales. Theorem~\ref{teo:identificaciones-H-F-W-entero} gives, for every $v\in\mathcal S'(\mathbb R^n)$,
\[
d_{s,p}\|v\|_{F^s_{p,2}(\mathbb R^n)}
\leq
\|v\|_{H^{s,p}(\mathbb R^n)}
\leq
D_{s,p}\|v\|_{F^s_{p,2}(\mathbb R^n)}.
\]
The constants depend only on $n,s,p$ and are therefore the same in all charts. Apply the inequalities to $\widetilde{(h_i u)\circ\phi_i^{-1}}$, raise to the power $p$, sum over $i$, and take the $p$th root. This proves \eqref{eq:equivalencia-local-H-F}; independence of the trivialization was proved earlier.

We now prove \eqref{eq:equivalencia-local-intrinseca-bessel}. First suppose $s>0$ and choose an integer $m>s$. At order zero, Proposition~\ref{prop:covariantes-localizacion-orden-entero} gives
\[
H^{0,p}(M)=L^p(M)=H_L^{0,p}(M).
\]
At order $m$, the same proposition and Theorem~\ref{teo:yoshida-bessel-covariante-geometria-acotada}, applied to the trivial rank-one bundle, give
\begin{equation}
\label{eq:extremo-entero-local-bessel-cap13}
H^{m,p}(M)=W^{m,p}(M)=H_L^{m,p}(M)
\end{equation}
and the comparisons in \eqref{eq:equivalencia-covariante-localizada-haz}, \eqref{eq:W-H-entero-geometria-acotada}, and \eqref{eq:yoshida-bessel-covariante} provide positive constants in both directions, independent of $u$. Set $\theta=\frac{s}{m}$. Proposition~\ref{prop:retraccion-interpolacion-dualidad-geometria-acotada} gives
\[
[H^{0,p}(M),H^{m,p}(M)]_\theta=H^{s,p}(M),
\]
while Theorem~\ref{teo:interpolacion-compleja-bessel-laplaciano-completo} gives
\[
[H_L^{0,p}(M),H_L^{m,p}(M)]_\theta=H_L^{s,p}(M).
\]
The identity on $C_c^\infty(M)$ is continuous between the two pairs at the endpoints by \eqref{eq:extremo-entero-local-bessel-cap13}. Theorem~\ref{teo:interpolacion-compleja-operadores}, applied to the identity and its inverse, proves equality and both estimates for $s>0$. For $s=0$, the assertion has already been obtained.

Finally, let $s<0$ and write $t=-s>0$. Applying the positive case with conjugate exponent $p'$ and taking duals, Proposition~\ref{prop:retraccion-interpolacion-dualidad-geometria-acotada} and the duality item of Theorem~\ref{teo:propiedades-escala-bessel-haces} give
\[
H^{-t,p}(M)
=
\bigl(H^{t,p'}(M)\bigr)'
=
\bigl(H_L^{t,p'}(M)\bigr)'
=
H_L^{-t,p}(M).
\]
In the real scalar case, both identifications use exactly the distributional pairing induced by $d\lambda_{\mathbf{g}}$. In the complex scalar case, the duality of Theorem~\ref{teo:propiedades-escala-bessel-haces} was written using the Hermitian product; composing the representing element with conjugation gives the bilinear distributional pairing of Proposition~\ref{prop:retraccion-interpolacion-dualidad-geometria-acotada}. The scalar semigroup and its potentials commute with conjugation, so this modification changes neither the space nor its norm. Thus, the equality holds in $\mathcal D'(M)$, not merely as an abstract Banach-space isomorphism. The norms of the isomorphisms and their duals produce the constants in \eqref{eq:equivalencia-local-intrinseca-bessel}.

The local comparisons and the passage to negative orders are justified by the preceding propositions.
\end{proof}

\begin{corollary}[Agreement of the constructions at integer orders]
\label{cor:sobolev-entero-cuatro-definiciones-geometria-acotada}
Let $(M,\mathbf{g})$ be a Riemannian manifold without boundary and with bounded geometry, $m\in\mathbb N_0$, and $1<p<\infty$. Then
\[
W^{m,p}(M)
=
H^{m,p}(M)
=
F^m_{p,2}(M)
=
H_L^{m,p}(M)
\]
and there exist constants $c_H,C_H,c_F,C_F,c_L,C_L>0$, independent of $u$, such that
\[
 c_H\|u\|_{W^{m,p}(M)}
 \leq\|u\|_{H^{m,p}(M)}
 \leq C_H\|u\|_{W^{m,p}(M)},
\]
\[
 c_F\|u\|_{W^{m,p}(M)}
 \leq\|u\|_{F^m_{p,2}(M)}
 \leq C_F\|u\|_{W^{m,p}(M)},
 \qquad
 c_L\|u\|_{W^{m,p}(M)}
 \leq\|u\|_{H_L^{m,p}(M)}
 \leq C_L\|u\|_{W^{m,p}(M)}.
\]

More generally, if $\mathbf{E}\to M$ is a bundle of bounded geometry, let $\mathbf{u}$ be a section of $\mathbf{E}$. Denote by $N_{m,p}(\mathbf{u})$ the localized component norm defined in \eqref{eq:equivalencia-covariante-localizada-haz}, and by $c_{m,p},C_{m,p},c^L_{m,p},C^L_{m,p}>0$ the constants in \eqref{eq:equivalencia-covariante-localizada-haz} and \eqref{eq:yoshida-bessel-covariante}. Then
\[
\begin{aligned}
c_{m,p}\|\mathbf u\|_{W^{m,p}(M,\mathbf E)}
&\leq N_{m,p}(\mathbf u)
 \leq C_{m,p}\|\mathbf u\|_{W^{m,p}(M,\mathbf E)},\\
c^L_{m,p}\|\mathbf u\|_{W^{m,p}(M,\mathbf E)}
&\leq\|\mathbf u\|_{H_{L_{\mathbf E}}^{m,p}(M,\mathbf E)}
 \leq C^L_{m,p}\|\mathbf u\|_{W^{m,p}(M,\mathbf E)}.
\end{aligned}
\]
The constants do not depend on the section.
\end{corollary}

\begin{proof}
In the scalar case, combine Theorem~\ref{teo:triebel-bessel-localizacion-geometria-acotada}, Proposition~\ref{prop:covariantes-localizacion-orden-entero}, and the Euclidean identification $H^{m,p}(\mathbb R^n)=F^m_{p,2}(\mathbb R^n)$. The vector-valued assertion follows directly by combining \eqref{eq:equivalencia-covariante-localizada-haz} and \eqref{eq:yoshida-bessel-covariante}.
\end{proof}

\begin{remark}[Uniformity and compactness]
\label{obs:uniformidad-no-compacidad-geometria-acotada}
The preceding norm comparisons allow local estimates to be summed with uniform constants. When the integrability exponent changes, global embeddings also use the inclusion of outer sums $\ell^p(I)\subseteq\ell^q(I)$ for $p\leq q$; comparison of differential indices does not replace this condition on a manifold of infinite volume. The precise ranges are formulated in Theorem~\ref{teo:encajes-escalas-geometria-acotada}.

Nor does a strict gain in index imply global compactness. Already on $\mathbb R^n$, translates of a nonzero smooth compactly supported function preserve its Sobolev norms and can have pairwise disjoint supports; their distances in $L^q(\mathbb R^n)$ remain equal to $2^{1/q}\|u\|_{L^q(\mathbb R^n)}>0$ for $q<\infty$. Proposition~\ref{prop:compacidad-local-encajes-geometria-acotada} obtains compactness after multiplication by a compactly supported cutoff and for families supported in a common compact set. A global conclusion on a noncompact manifold also requires control of mass escaping to infinity, as explained in Remark~\ref{obs:perdida-compacidad-global-geometria-acotada}.
\end{remark}

\section{Laplacian spectrum in bounded geometry}
\label{sec:espectro-laplaciano-geometria-acotada}

On a noncompact manifold, local elliptic regularity and uniform estimates
do not provide a globally compact inclusion. We therefore distinguish
eigenvalues, discrete spectrum, and essential spectrum.
The Laplacian realization remains the nonnegative self-adjoint one
constructed in Chapter~\ref{cap:sobolev-fraccionario-haces-nucleo-calor}.
The first assertions concern self-adjoint operators; we then apply
them to the Laplacian and its perturbations.

\begin{definition}[Essential spectrum and spectral gap]
\label{def:espectro-esencial-brecha-cap24}\glsadd{espectro-esencial}
\index{essential spectrum}\index{spectral gap}
Let $A$ be a self-adjoint operator on a complex Hilbert
space $H$. Its \textit{essential spectrum} is
\[
 \sigma_{\mathrm e}(A):=
 \{\lambda\in\mathbb R\mid
   A-\lambda I:\mathcal D(A)\longrightarrow H
   \text{ is not Fredholm}\},
\]
where $\mathcal D(A)$ carries the graph norm.
The \textit{discrete spectrum} consists of isolated eigenvalues of
finite multiplicity. If $E_A$ is the spectral measure of $A$, write
$\mu_u(B):=\langle E_A(B)u,u\rangle_H$.

We say that $A$ has a \textit{spectral gap above its kernel} if
there exists $\gamma>0$ such that
\[
 \langle Au,u\rangle_H\geq\gamma\|u\|_H^2,
 \qquad u\in\mathcal D(A)\cap(\ker A)^\perp.
\]
Write $A_\perp$ for its restriction to the reducing subspace
$H_\perp:=(\ker A)^\perp$, with domain $\mathcal D(A)\cap H_\perp$.
\end{definition}

\begin{theorem}[Spectral projections and Weyl criterion]
\label{teo:espectro-esencial-proyecciones-weyl}
For $A$ as above and $\lambda\in\mathbb R$, the following are equivalent:
\begin{enumerate}[label=(\alph*)]
 \item $\lambda\in\sigma_{\mathrm e}(A)$;
 \item $E_A((\lambda-\varepsilon,\lambda+\varepsilon))$ has
 infinite-dimensional range for every $\varepsilon>0$;
 \item there exist $u_j\in\mathcal D(A)$ such that
 \[
  \|u_j\|_H=1,\qquad u_j\rightharpoonup0,\qquad
  \|(A-\lambda I)u_j\|_H\longrightarrow0.
 \]
\end{enumerate}
The sequence in (c) is called a \textit{singular Weyl sequence}.
Moreover,
\[
 \sigma(A)=\sigma_{\mathrm{disc}}(A)\,\dot\cup\,\sigma_{\mathrm e}(A).
\]
An isolated point of the spectrum of a self-adjoint operator is an
eigenvalue, although its multiplicity may be infinite.
\end{theorem}

\begin{proof}
Set $T:=A-\lambda I$ and
$P_\varepsilon:=E_A((\lambda-\varepsilon,\lambda+\varepsilon))$.
Spectral calculus gives, for $u\in\mathcal D(A)$,
\[
 \|Tu\|_H^2=\int_{\mathbb R}|t-\lambda|^2\,d\mu_u(t),\qquad
 \|T(I-P_\varepsilon)u\|_H\geq
                  \varepsilon\|(I-P_\varepsilon)u\|_H.
\]
The range of $P_\varepsilon$ is contained in $\mathcal D(A)$ because
the interval is bounded. On its orthogonal complement, the
multiplier $(t-\lambda)^{-1}$ defines a continuous inverse of
$T$ into its graph-norm domain.

If $P_\varepsilon$ has finite rank, the restriction of $T$ to its
range is an operator between finite-dimensional spaces, and the
restriction to the complement is invertible. Thus $T$ has finite
kernel and cokernel and closed range, so it is Fredholm.
Conversely, if $T$ is Fredholm, its kernel $N$ is finite dimensional,
and the map
\[
 T:\mathcal D(A)\cap N^\perp\longrightarrow\operatorname{Ran}T
\]
is a Banach isomorphism. Indeed, its range is closed, and injectivity
and surjectivity follow from
$\operatorname{Ran}T=(\ker T)^\perp$ by self-adjointness.
Continuity of the inverse gives $\|u\|_H\leq C\|Tu\|_H$ on that domain.
Choose $0<\varepsilon<C^{-1}$.
For $u\in\operatorname{Ran}P_\varepsilon\cap N^\perp$, we also have
$\|Tu\|_H\leq\varepsilon\|u\|_H$, so $u=0$.
Since $P_\varepsilon$ contains $N$ and commutes with its projection,
its range equals $N$ and is finite dimensional.
This proves the equivalence of (a) and (b).

If (b) holds, choose inductively a unit vector $u_j$ in
$\operatorname{Ran}P_{1/j}$ orthogonal to $u_1,\ldots,u_{j-1}$.
This is possible because that range is infinite dimensional.
The sequence is orthonormal, hence weakly null by Bessel, and
$\|Tu_j\|_H\leq1/j$. This proves (c).
If (c) holds and some $P_\varepsilon$ has finite rank, weak
convergence implies $\|P_\varepsilon u_j\|_H\to0$. Then
\[
 \|Tu_j\|_H\geq\varepsilon\|(I-P_\varepsilon)u_j\|_H
 \longrightarrow\varepsilon,
\]
a contradiction. The three conditions are equivalent.

To identify the spectrum, a real number $\lambda$ belongs to
$\sigma(A)$ if and only if $P_\varepsilon\neq0$ for every
$\varepsilon>0$. If a projection is zero, the multiplier
$(t-\lambda)^{-1}$ is bounded on the spectral support and gives
the resolvent. Conversely, if $(A-\lambda I)^{-1}$ has norm $C$
and $P_\varepsilon\neq0$ with $C\varepsilon<1$, a unit vector in
its range contradicts $\|u\|_H\leq C\|Tu\|_H$.
If $\lambda\in\sigma(A)$ and $P_\varepsilon$ has finite rank,
the self-adjoint operator on that range diagonalizes. Its spectrum
is finite, and the spectrum on the complement lies outside the
interval. Hence $\lambda$ is an isolated eigenvalue of finite
multiplicity. Conversely, for such an eigenvalue, a projection on
a sufficiently small interval is the projection onto its
eigenspace and has finite rank.
If $\lambda$ is isolated, that projection equals
$E_A(\{\lambda\})$ and is nonzero; its range is
$\ker(A-\lambda I)$. This also proves the final assertion. The essential spectrum is
closed: if a projection on an interval about $\lambda$ has finite rank,
the same holds on smaller intervals about points sufficiently close
to $\lambda$.
\end{proof}

\begin{remark}[Essential and continuous spectrum]
\label{obs:esencial-no-equivale-continuo}
The essential spectrum can contain eigenvalues: the zero operator
on an infinite-dimensional Hilbert space has zero both as an
eigenvalue and as its sole essential spectral point.
For a self-adjoint operator, if $\lambda\in\sigma(A)$ is not
an eigenvalue, $A-\lambda I$ is injective, has dense range, and
is not surjective. Density follows from
$\operatorname{Ran}(A-\lambda I)^\perp=\ker(A-\lambda I)$.
If it were surjective, its inverse would be bounded by the closed
graph theorem and $\lambda$ would belong to the resolvent set.
These are the continuous spectral points. The Weyl condition
also detects infinite-multiplicity eigenvalues and accumulations
of discrete spectrum.
\end{remark}

\begin{proposition}[Rellich criterion for the spectrum]
\label{prop:rellich-resolvente-espectro-esencial}
For a nonnegative self-adjoint operator $A$, the following are equivalent:
\begin{enumerate}[label=(\alph*)]
 \item $\sigma_{\mathrm e}(A)=\varnothing$;
 \item $(I+A)^{-1}$ is compact on $H$;
 \item the inclusion $\mathcal D(A^{1/2})\hookrightarrow H$ is compact
 for the norm
 $\|u\|_V^2=\|u\|_H^2+\|A^{1/2}u\|_H^2$.
\end{enumerate}
\end{proposition}

\begin{proof}
If the essential spectrum is empty, each point of
$\sigma(A)\cap[0,N]$ has an interval whose spectral projection
has finite rank. That spectral set is compact and is covered
by finitely many such intervals.
Thus $P_N:=E_A([0,N])$ has finite rank: its range is contained
in the sum of those ranges. Spectral calculus gives
\[
 \|(I+A)^{-1}-(I+A)^{-1}P_N\|\leq(1+N)^{-1}.
\]
The resolvent is a norm limit of finite-rank operators.
Conversely, if the resolvent is compact and some $P_N$ has
infinite rank, take an orthonormal sequence $(u_j)$ in its range.
On that range,
\[
 \|(I+A)^{-1}(u_j-u_k)\|_H
 \geq(1+N)^{-1}\|u_j-u_k\|_H
 =\frac{\sqrt2}{1+N}\quad(j\neq k),
\]
contradicting compactness. All $P_N$ have finite rank, and the
projection criterion excludes essential spectrum.

The map $U:V\to H$, $Uu=(I+A)^{1/2}u$, is unitary for the stated
norm, with inverse $S=(I+A)^{-1/2}$. The inclusion $J:V\to H$
satisfies $J=SU$ and $S=JU^{-1}$, so $J$ is compact if and only
if $S$ is. If $S$ is compact, then $S^2=(I+A)^{-1}$ is compact.
If $S^2$ is compact, its compact self-adjoint spectral theorem
shows that its positive square root is a norm limit of finite-rank
operators: take square roots of its eigenvalues, which tend to zero.
That root is $S$ by functional calculus. This proves equivalence
with (c).
\end{proof}

\begin{proposition}[Gap, closed range, and inverse]
\label{prop:brecha-espectral-rango-inversa}
Let $A$ be nonnegative and self-adjoint. If $H_\perp\neq\{0\}$, the following conditions are equivalent:
\begin{enumerate}[label=(\alph*)]
 \item there is a gap $\gamma>0$ above the kernel;
 \item $\sigma(A_\perp)\subset[\gamma,\infty)$ for some $\gamma>0$;
 \item $\operatorname{Ran}A$ is closed.
\end{enumerate}
In that case, $\operatorname{Ran}A=H_\perp$, the inverse
$G:=A_\perp^{-1}:H_\perp\to\mathcal D(A)\cap H_\perp$ is continuous,
and $\|G\|_{\mathcal L(H_\perp)}\leq\gamma^{-1}$.
The gap inequality extends to $\mathcal D(A^{1/2})\cap H_\perp$:
\[
 \|A^{1/2}u\|_H^2\geq\gamma\|u\|_H^2.
\]
In particular, $\inf\sigma_{\mathrm e}(A_\perp)>0$ implies these
conditions, with $\inf\varnothing=+\infty$.
If $H_\perp=\{0\}$, the range and inverse assertions hold with
the zero inverse.
\end{proposition}

\begin{proof}
A gap on $\mathcal D(A)\cap H_\perp$ forces
$E_{A_\perp}([0,\gamma-\varepsilon])=0$ for
$0<\varepsilon<\gamma$: a nonzero vector in that range belongs
to the domain and has Rayleigh quotient at most
$\gamma-\varepsilon$. The spectral support criterion in the
preceding proof gives (b). Conversely, (b) and
$\langle Au,u\rangle=\int t\,d\mu_u(t)$ give (a) and its
extension to the square-root domain.

Under (b), the multiplier $t^{-1}$ on $[\gamma,\infty)$ provides
$G$, with $\|Gf\|_H\leq\gamma^{-1}\|f\|_H$ and $AGf=f$.
Its range is in $\mathcal D(A)$ because $t\,t^{-1}=1$ on that
support. Hence $\operatorname{Ran}A=H_\perp$ is closed.
If the range is closed, self-adjointness gives
$\operatorname{Ran}A=H_\perp$. The bounded inverse theorem for
the graph norm provides $C>0$ with $\|u\|_H\leq C\|Au\|_H$
on $\mathcal D(A)\cap H_\perp$. Projections onto $[0,\varepsilon]$
must then vanish for $\varepsilon<C^{-1}$, by the same norm
comparison used in the Weyl criterion. This proves (b).

Finally, if $\inf\sigma_{\mathrm e}(A_\perp)>0$, zero is not
essential. If it belonged to $\sigma(A_\perp)$, it would be an
isolated eigenvalue of finite multiplicity by
Theorem~\ref{teo:espectro-esencial-proyecciones-weyl}.
This contradicts $\ker A_\perp=\{0\}$. Thus zero belongs to
the resolvent. Since the spectrum is closed and nonnegative,
there exists $\gamma>0$ with
$\sigma(A_\perp)\subset[\gamma,\infty)$.
\end{proof}

\begin{proposition}[Rayleigh and minimax below essential spectrum]
\label{prop:minimax-umbral-espectro-esencial}
Let $A$ be a nonnegative self-adjoint operator on an
infinite-dimensional Hilbert space. Let $\mathfrak q$ be its form,
with domain $V=\mathcal D(A^{1/2})$, and set
\[
 \alpha_k:=\inf_{\substack{W\subset V\\\dim W=k}}
        \ \max_{u\in W\setminus\{0\}}
             \frac{\mathfrak q(u,u)}{\|u\|_H^2},
 \qquad b:=\inf\sigma_{\mathrm e}(A),\quad\inf\varnothing:=+\infty.
\]
If there are at least $k$ eigenvalues below $b$, counted with
multiplicity, $\alpha_k$ is the $k$th such eigenvalue.
If there are fewer than $k$, then $b<\infty$ and $\alpha_k=b$.
Moreover,
\begin{equation}
 \inf\sigma(A)=\inf_{u\in V\setminus\{0\}}
                   \frac{\mathfrak q(u,u)}{\|u\|_H^2}.
 \label{eq:fondo-espectro-rayleigh}
\end{equation}
The last infimum is attained if and only if the spectral bottom
is an eigenvalue.
\end{proposition}

\begin{proof}
For $c\in\mathbb R$, write $P_c:=E_A([0,c))$, with $P_c=0$ if
$c\leq0$. If $\operatorname{rank}P_c<k$, every $k$-dimensional
subspace $W$ contains a nonzero vector orthogonal to that range.
For $c>0$, its spectral measure is supported in $[c,\infty)$,
and $\mathfrak q(u,u)\geq c\|u\|_H^2$. For $c\leq0$ this
inequality follows from nonnegativity. Hence $\alpha_k\geq c$.
If $\operatorname{rank}P_c\geq k$, choose a $k$-dimensional
subspace $W$ inside its range. There $A$ is bounded and
$\mathfrak q(u,u)\leq c\|u\|_H^2$, so $\alpha_k\leq c$.
Consequently,
\[
 \alpha_k=\inf\{c\in\mathbb R\mid\operatorname{rank}P_c\geq k\}.
\]
The set is nonempty because $P_c\to I$ strongly as
$c\to+\infty$; if all its ranges had dimension less than $k$,
no set of $k$ independent vectors could be approximated by them.
More explicitly, choose orthonormal $e_1,\ldots,e_k$ and $c$ such that
$\|(I-P_c)e_i\|_H<1/(2\sqrt{k})$ for every $i$. For
$\sum_i|a_i|^2=1$,
\[
 \left\|\sum_{i=1}^ka_iP_ce_i\right\|_H
 \geq1-\sum_{i=1}^k|a_i|\|(I-P_c)e_i\|_H>\frac12.
\]
Their images are independent, so the range has dimension at least $k$.

For $0\leq c<b$, the projection $E_A([0,c])$ has finite rank.
Indeed, the compact set $\sigma(A)\cap[0,c]$ is covered by
finitely many intervals with finite-rank projections, as in the
Rellich criterion. All its spectral points are therefore
isolated eigenvalues of finite multiplicity.
The rank formula identifies $\alpha_k$ with the $k$th eigenvalue
when it lies below $b$. If there are fewer than $k$ eigenvalues
below $b$, no $P_c$ with $c<b$ has rank at least $k$.
The case $b=+\infty$ would contradict the strong convergence
above. If $b<\infty$, closedness of the essential spectrum
gives $b\in\sigma_{\mathrm e}(A)$; every
$P_{b+\varepsilon}$ has infinite rank by the Weyl criterion.
Thus the rank formula gives $\alpha_k=b$.

For $k=1$, the same support criterion gives
\eqref{eq:fondo-espectro-rayleigh}.
If $a:=\inf\sigma(A)$, equality in the quotient requires
$\int(t-a)\,d\mu_u(t)=0$. Since the integrand is nonnegative,
the measure is supported on $\{a\}$, and
$u\in\ker(A-aI)$. An eigenvector of $a$ gives equality conversely.
\end{proof}

\begin{proposition}[Volume, kernel, and local compactness]
\label{prop:laplaciano-abierto-volumen-nucleo-resolvente}
Let $(M,\mathbf g)$ be connected, noncompact, without boundary,
and of bounded geometry of all orders.
Then $\lambda_{\mathbf g}(M)=\infty$, and the scalar Laplacian
$L=-\Delta_{\mathbf g}$ satisfies:
\begin{enumerate}[label=(\alph*)]
 \item $\ker_{L^2}L=\{0\}$;
 \item $\chi(I+L)^{-1}$ is compact on $L^2(M)$ for
       $\chi\in C_c^\infty(M)$;
 \item $(I+L)^{-1}$ is not compact and
       $\sigma_{\mathrm e}(L)\neq\varnothing$.
\end{enumerate}
The volume and kernel conclusions also hold under bounded geometry
of order zero: positive injectivity radius and uniformly bounded curvature.
For a bundle with compatible connection, the $L^2$ kernel of its
connection Laplacian is likewise zero on a connected base
of infinite volume.
\end{proposition}

\begin{proof}
Choose $0<\rho<\operatorname{inj}(M)/4$ within the uniform normal
radius. The injectivity-radius lower bound extends each unit-speed
geodesic by intervals of length $\rho$ from any of its points.
A finite maximal time is impossible; the manifold is geodesically
complete, and Hopf--Rinow applies.
Lemma~\ref{lema: volumen bolas geometria acotada} provides
$v_->0$ with $\lambda_{\mathbf g}(B_{\mathbf g}(p,\rho))\geq v_-$
for every $p$. Choose centers $p_j$ inductively with
$d_{\mathbf g}(p_j,p_i)>3\rho$ for $i<j$.
After finitely many choices, the union of the closed balls of
radius $3\rho$ is compact by Hopf--Rinow and cannot cover the
noncompact base. Thus the next center exists. The balls of radius
$\rho$ are disjoint, and total volume is infinite.

For that conclusion, bounded geometry of order zero suffices.
Indeed, if $|\mathbf R|\leq C_0$, also choose $C_0\rho^2\leq1$.
In a parallel frame along a unit-speed geodesic, the normal Jacobi
matrix satisfies
\[
 J(t)=tI-\int_0^t(t-s)\mathcal R(s)J(s)\,ds,
 \qquad\|\mathcal R(s)\|\leq C_0.
\]
For $h_t:=\sup_{0<s\leq t}\|J(s)\|/s$, this gives
$h_t\leq1+C_0t^2h_t/6$, hence $h_t\leq2$ and
$\|J(t)/t-I\|\leq C_0t^2/3\leq1/3$ for $t\leq\rho$.
The singular values of $J(t)/t$ are at least $1/2$;
the radial factor is one. The volume density in normal coordinates
is therefore at least $2^{-(n-1)}$.
Thus $\lambda_{\mathbf g}(B(p,\rho))\geq
2^{-(n-1)}\omega_n\rho^n$ uniformly in $p$, where $\omega_n$ is the
Euclidean unit-ball volume. For $n=1$, the radial density is one.
This justifies the same ball construction under order-zero assumptions.

If $u\in\ker L$, the form identity gives $\|du\|_{L^2}=0$.
Local elliptic regularity makes $u$ smooth, hence constant by
connectedness. Integrability and infinite volume force $u=0$.
For a connection Laplacian, the same argument gives
$\nabla\mathbf u=0$; metric compatibility makes $|\mathbf u|$
constant, and that constant must be zero.

If $f_j$ is bounded in $L^2$ and $u_j=(I+L)^{-1}f_j$, then
$u_j$ and $Lu_j=f_j-u_j$ are bounded in $L^2$.
Choose a relatively compact smooth domain $D$ containing the
support of $\chi$, and another $D'$ with
$\overline D\subset D'$ and compact $\overline{D'}$.
Interior elliptic estimates bound $\|u_j\|_{H^2(D)}$ by the
norms of $u_j$ and $Lu_j$ on $D'$.
Multiplication by $\chi$ and Rellich on $D$ yield an $L^2(M)$
convergent subsequence of $\chi u_j$.
This proves (b) without global compactness.

For (c), take a nonzero
$\eta\in C_c^\infty(B_{\mathrm{euc}}(0,\rho/2))$ and transport
it to the normal charts centered at $p_j$, extending by zero.
Call the resulting functions $w_j$. Their supports lie strictly
inside the charts. Uniform bounds for the metric and its inverse,
the volume density, and first derivatives give
\[
 0<c\leq\|w_j\|_{L^2}\leq C,\qquad
 \|Lw_j\|_{L^2}\leq C.
\]
The last estimate follows from the local formula for $-\Delta$:
only derivatives of $\eta$ up to order two and uniformly bounded
coefficients occur.
The vectors $v_j:=w_j/\|w_j\|_{L^2}$ are orthonormal by disjoint
supports and uniformly bounded in graph norm.
If $(I+L)^{-1}$ were compact, the bounded sequence $(I+L)v_j$
would yield a convergent subsequence of
$v_j=(I+L)^{-1}(I+L)v_j$, a contradiction.
The Rellich criterion gives nonempty essential spectrum.
\end{proof}

\begin{corollary}[Spectral bottom from compactly supported tests]
\label{cor:fondo-laplaciano-testigos-compactos}
On a connected complete manifold without boundary,
\[
 \inf\sigma(-\Delta_{\mathbf g})
 =\inf_{u\in C_c^\infty(M)\setminus\{0\}}
       \frac{\|du\|_{L^2}^2}{\|u\|_{L^2}^2}.
\]
The same formula holds for the connection Laplacian and compactly
supported smooth sections.
On a connected noncompact base of bounded geometry, zero is
not an $L^2$ eigenvalue of the connection Laplacian; this
does not exclude zero from its spectrum.
\end{corollary}

\begin{proof}
The form domain is $W^{1,2}$, and compactly supported smooth
functions or sections are dense for its norm, by completeness
and the cutoffs constructed in the heat chapter.
For a nonzero form-domain element, approximation in that norm
gives denominators converging to its positive $L^2$ norm squared
and numerators converging to its energy.
Density and \eqref{eq:fondo-espectro-rayleigh} prove equality of
the infima. The kernel assertion was proved above.
\end{proof}

\begin{proposition}[Compact perturbations in graph norm]
\label{prop:invariancia-esencial-perturbacion-grafica}
Let $A$ and $B$ be self-adjoint operators on $H$ with the same
domain $D$. If $B-A:D\to H$ is compact for the graph norm of $A$,
then $\sigma_{\mathrm e}(A)=\sigma_{\mathrm e}(B)$.
In particular, on a manifold without boundary of bounded
geometry, if $V\in L^\infty(M)$ is real and tends uniformly to
zero outside compact sets, then
\[
 \sigma_{\mathrm e}(-\Delta_{\mathbf g}+V)
             =\sigma_{\mathrm e}(-\Delta_{\mathbf g}).
\]
The perturbed operator has the same domain as $L$.
\end{proposition}

\begin{proof}
The graph norms of $A$ and $B$ are equivalent.
Indeed, the identity between the two Banach spaces has closed
graph: convergence in either norm implies convergence to
the same vector in $H$. The closed graph theorem in both
directions gives equivalence.
For $\lambda\in\mathbb R$, the maps
$A-\lambda I,B-\lambda I:D\to H$ differ by a compact operator.
Stability of Fredholm operators under compact perturbations,
applied in both directions, says that one is Fredholm if and
only if the other is. The essential spectrum definition gives
the conclusion. For self-adjoint operators bounded below, the
same proof applies after a common shift making them nonnegative.

Multiplication by $V$ is bounded and self-adjoint on $L^2(M)$.
By the adjoint definition,
$(L+V)^*=L^*+V$ with domain
$\mathcal D(L^*)=\mathcal D(L)$:
subtracting the bounded functional associated with $V$ does
not change the condition defining that domain.
Hence $L+V$ is self-adjoint and
$L+V\geq-\|V\|_\infty I$.

Let $K_j$ be an exhaustion by compact sets with
$\sup_{M\setminus K_j}|V|\to0$, and choose
$\chi_j\in C_c^\infty(M)$ equal to one near $K_j$, with
$0\leq\chi_j\leq1$. Multiplication by $\chi_j$ is compact
from $\mathcal D(L)$ with graph norm to $L^2(M)$:
apply local elliptic estimates and Rellich, as in (b) above.
Composition with bounded multiplication by $V$ makes
$\chi_jV$ compact. Moreover,
\[
 \|(1-\chi_j)Vu\|_{L^2}
 \leq\sup_{M\setminus K_j}|V|\,\|u\|_{L^2}
 \leq\sup_{M\setminus K_j}|V|\,\|u\|_{\mathcal D(L)}.
\]
Thus $V:\mathcal D(L)\to L^2(M)$ is a norm limit of compact
operators. The first part applies.
\end{proof}

\begin{example}[Euclidean Laplacian]
\label{ej:espectro-esencial-laplaciano-euclidiano}
For $L=-\Delta$ on $\mathbb R^n$, $n\geq1$,
\[
 \sigma(L)=\sigma_{\mathrm e}(L)=[0,\infty),\qquad
 \sigma_{\mathrm p}(L)=\varnothing.
\]
The unitary Fourier transform identifies $L$ with multiplication
by $|\xi|^2$. If $Lu=\lambda u$, its transform $\widehat u$
is supported on $\{|\xi|^2=\lambda\}$, a null set for
$\lambda\geq0$ and an empty set for $\lambda<0$.
An $L^2$ function with that support is zero.

For $\lambda=|\xi_0|^2\geq0$, fix
$\chi\in C_c^\infty(B_{\mathrm{euc}}(0,1))$ with $\|\chi\|_{L^2}=1$,
$R_j\to\infty$, and centers $a_j$ such that
$B(a_j,R_j)$ are disjoint. Choose them inductively outside
the finite union of the previously chosen enlarged balls.
Set
\[
 u_j(x):=R_j^{-n/2}\chi((x-a_j)/R_j)e^{i\langle\xi_0,x\rangle}.
\]
These vectors are orthonormal, and
\[
 \|(L-\lambda I)u_j\|_{L^2}
 \leq R_j^{-2}\|\Delta\chi\|_{L^2}
          +2|\xi_0|R_j^{-1}\|d\chi\|_{L^2}\longrightarrow0.
\]
The Weyl criterion places every $\lambda\geq0$ in the essential
spectrum. Nonnegativity excludes negative values, and
self-adjointness excludes nonreal spectrum.
This proves the stated equalities.
In particular, the spectral bottom zero is not an eigenvalue,
and the Rayleigh infimum is not attained.
\end{example}

\begin{example}[Positive spectral bottom in bounded geometry]
\label{ej:espectro-esencial-hiperbolico}
Consider the horospherical model of hyperbolic space, $n\geq2$,
\[
 M=\mathbb R\times\mathbb R^{n-1},\qquad
 \mathbf g=ds^2+e^{2s}|dx|^2,\qquad \rho:=\frac{n-1}{2}.
\]
It is complete. If $(s_j,x_j)$ is Cauchy, the bound
$|ds|\leq|\cdot|_{\mathbf g}$ makes $s_j$ Cauchy; choose $B$ with
$|s_j|\leq B$. For sufficiently large $j,k$ and small
$\varepsilon>0$, choose a curve joining the points with length less than
$d_{\mathbf g}((s_j,x_j),(s_k,x_k))+\varepsilon<1$.
Along the entire curve, $|s|\leq B+1$, so
\[
 |x_j-x_k|\leq e^{B+1}
       \bigl(d_{\mathbf g}((s_j,x_j),(s_k,x_k))+\varepsilon\bigr).
\]
Letting $\varepsilon\downarrow0$ makes $x_j$ Cauchy.
Both coordinates converge, and local comparison with the Euclidean
metric gives convergence in Riemannian distance.
The maps $(s,x)\mapsto(s+a,e^{-a}x+b)$ are isometries acting
transitively. Thus the injectivity radius is a positive constant,
and the norms of all covariant curvature derivatives are finite
constants: the geometry is bounded.

The volume element and positive Laplacian are
\[
 d\lambda_{\mathbf g}=e^{2\rho s}\,ds\,dx,\qquad
 L=-\partial_s^2-2\rho\partial_s-e^{-2s}\Delta_x.
\]
The field $X=\partial_s$ has norm one and
$\operatorname{div}X=2\rho$.
For $u\in C_c^\infty(M)$, integration by parts and
Cauchy--Schwarz give
\[
 2\rho\|u\|_{L^2}^2
 =-2\operatorname{Re}\int_M\overline u\,Xu\,d\lambda_{\mathbf g}
 \leq2\|u\|_{L^2}\|du\|_{L^2}.
\]
Thus $\sigma(L)\subset[\rho^2,\infty)$ by the Rayleigh formula.

The unitary transformation $u=e^{-\rho s}v$ into $L^2(ds\,dx)$
converts $L$ to $-\partial_s^2+\rho^2-e^{-2s}\Delta_x$.
Choose $\chi\in C_c^\infty((0,1))$ and
$\eta\in C_c^\infty(\mathbb R^{n-1})$, both of $L^2$ norm one.
For $\lambda=\rho^2+k^2$, $k\geq0$, define
\[
 v_j(s,x):=j^{-1/2}
      \chi((s-j^3)/j)e^{iks}\eta(x),\qquad
 u_j:=e^{-\rho s}v_j.
\]
The supports lie in disjoint bands $j^3<s<j^3+j$;
$u_j$ are orthonormal in $L^2(M)$.
Direct differentiation gives
\[
 \|(L-\lambda I)u_j\|_{L^2(M)}
 \leq j^{-2}\|\chi''\|_{L^2}
       +2k j^{-1}\|\chi'\|_{L^2}
       +e^{-2j^3}\|\Delta_x\eta\|_{L^2}\longrightarrow0.
\]
The Weyl criterion proves
$\sigma(L)=\sigma_{\mathrm e}(L)=[\rho^2,\infty)$.
Hence bounded geometry does not force zero into the spectrum;
here the Laplacian range is all of $L^2(M)$, and its inverse has
norm at most $\rho^{-2}$.
\end{example}

\begin{remark}[Connection and Hodge Laplacians]
\label{obs:espectro-bochner-hodge-geometria-acotada}
The projection, Weyl, minimax, and gap criteria concern self-adjoint
operators and also apply to the Hodge Laplacian
$\Delta_q=d\delta+\delta d$ with its self-adjoint realization.
Its form is $\|d\omega\|_{L^2}^2+\|\delta\omega\|_{L^2}^2$.
In positive degrees, this is not the connection energy:
the Weitzenb\"ock identity includes the curvature endomorphism.
Consequently, the argument that eliminates the connection
Laplacian kernel through the constant norm of a parallel section
does not in general eliminate the kernel of $\Delta_q$.
Fluid applications use $\Delta_1$, whereas gauge applications use
the covariant Laplacian on sections of $\operatorname{ad}\mathbf P$.
In both cases, the essential spectrum hypothesis concerns the
operator and kernel complement specified in the statement.
\end{remark}

Consultation sources for this section:
\cite{Eichhorn2007,KrupkaSaunders2008,Borthwick2020Spectral}.

\chapter{Trace theorems in bounded geometry}
\label{cap:trazas-geometria-acotada}

A smooth function can be restricted to a submanifold without difficulty. For a Sobolev function, however, its values on a set of measure zero are not determined by its equivalence class. The trace resolves this problem when there is sufficient transverse regularity: it is defined by continuity, identifies the space receiving the restriction, and allows us to ask whether every datum in that space admits an extension.

The threshold comes from the Euclidean model. To transfer it to a noncompact manifold, we must control adapted coordinates, frames, and cutoff functions, and justify the sums reconstructing the trace and its right inverse. Fermi charts keep the submanifold visible; in the presence of a boundary, the collar also allows tangential derivatives to be separated from normal jets.

We begin with these geometric constructions and their uniform estimates. We then treat traces of functions and sections, comparison with covariant norms, and Sobolev embeddings. The exposition distinguishes restriction spaces on the half-space from extensions by zero and preserves the order of summation defining each function-space scale; see \cite{traces,Schneider2021,Amann2025FunctionSpaces}.

\begin{semblanzaHistorica}{The boundary as analytic data}
Restricting a smooth function to a submanifold seems elementary; for a Sobolev class, however, it is possible only when there is sufficient transverse regularity. The trace theorem precisely identifies what information survives at the boundary and which space receives it. On noncompact manifolds, the problem also requires collars, Fermi charts, and extensions to be constructed with uniform constants. This combination of adapted geometry and fractional regularity is one of the contemporary forms of global analysis~\cite{traces,VanSchaftingenWinter2025}.
\end{semblanzaHistorica}

\section{The Euclidean model and its compatibility properties}

Let \(n>k\geq0\) and \(c:=n-k\). Identify \(\mathbb R^k\) with \(\{0\}^c\times\mathbb R^k\subseteq\mathbb R^c\times\mathbb R^k\), and, for a smooth function \(u\), write
\[
\operatorname{tr}_{c}u(x'):=u(0,x'),
\qquad x'\in\mathbb R^k.
\]
When $k=0$, $\mathbb R^0$ is a point with counting measure; its $H^{t,p}(\mathbb R^0)$, $F^t_{p,q}(\mathbb R^0)$, and $B^t_{p,q}(\mathbb R^0)$ spaces are identified with the scalar field, equipped with its absolute-value norm. This convention includes the boundary of a one-dimensional manifold. The Euclidean result proved in Theorem~\ref{teo:traza-fraccionaria-euclidea-hiperplano} will be used in the following form.

\begin{theorem}[Euclidean trace in arbitrary codimension]
\label{teo:traza-euclidiana-codimension-c}
Let \(1<p<\infty\), \(s>\frac{c}{p}\), and \(c=n-k\). Restriction of Schwartz functions extends uniquely to a continuous surjective linear operator
\[
\operatorname{tr}_{c}\colon H^{s,p}(\mathbb R^c\times\mathbb R^k)
\longrightarrow
B^{s-\frac{c}{p}}_{p,p}(\mathbb R^k).
\]
There exists a continuous linear operator
\[
\operatorname{ex}_{c}\colon B^{s-\frac{c}{p}}_{p,p}(\mathbb R^k)
\longrightarrow
H^{s,p}(\mathbb R^n)
\]
such that \(\operatorname{tr}_{c}\operatorname{ex}_{c}=\operatorname{id}\). The operators can be chosen to be continuous also on the scales
\[
\operatorname{tr}_{c}\colon F^s_{p,q}(\mathbb R^n)\longrightarrow
B^{s-\frac{c}{p}}_{p,p}(\mathbb R^k),
\qquad
\operatorname{tr}_{c}\colon B^s_{p,q}(\mathbb R^n)\longrightarrow
B^{s-\frac{c}{p}}_{p,q}(\mathbb R^k),
\]
for \(1\leq q\leq\infty\), with corresponding right extension operators.
\end{theorem}

\begin{proof}
For the last step when $k=0$, we need evaluation at a point of $\mathbb R$. If $A^s_{p,q}$ denotes $F^s_{p,q}$ or $B^s_{p,q}$, Theorem~\ref{teo:desigualdades-bernstein} in dimension one gives
\[
 \sum_{j\in\mathbb N_0}\|\Delta_j u\|_\infty
 \leq C\sum_{j\in\mathbb N_0}2^{j/p}\|\Delta_j u\|_{L^{p}(\mathbb R)}
 \leq C\left(\sum_{j\in\mathbb N_0}2^{-j(s-1/p)}\right)
            \|u\|_{A^s_{p,q}},\qquad s>1/p.
\]
In both scales, each $2^{js}\|\Delta_j u\|_{L^{p}(\mathbb R)}$ is bounded by the space norm. The series converges uniformly and defines a continuous representative whose evaluation at zero is a continuous trace. A fixed function $\varphi\in\mathcal S(\mathbb R)$ with $\varphi(0)=1$ gives the coretraction $a\mapsto a\varphi$, common to all orders. This also proves the Bessel case through $H^{s,p}(\mathbb R^n)=F^s_{p,2}(\mathbb R^n)$.

Write a point of \(\mathbb R^n\) as \((z_1,\dots,z_c,x')\), with \(x'\in\mathbb R^k\), and denote by \(\tau_r\) the restriction \(z_r=0\) after eliminating the variables \(z_1,\dots,z_{r-1}\). Theorem~\ref{teo:traza-fraccionaria-euclidea-hiperplano} and its coretraction give, for \(r\in\{1,\dots,c\}\),
\[
\tau_r\colon
F^{s-\frac{r-1}{p}}_{p,p}(\mathbb R^{n-r+1})
\longrightarrow
B^{s-\frac{r}{p}}_{p,p}(\mathbb R^{n-r})
=F^{s-\frac{r}{p}}_{p,p}(\mathbb R^{n-r}).
\]
The equality in the last member is \(B^\sigma_{p,p}(\mathbb R^d)=F^\sigma_{p,p}(\mathbb R^d)\). The condition \(s>\frac{c}{p}\) implies
\[
s-\frac{r-1}{p}>\frac{1}{p},
\qquad r\in\{1,\dots,c\},
\]
so every operator is defined. Since \(H^{s,p}(\mathbb R^n)=F^s_{p,2}(\mathbb R^n)\), the first restriction is applied with second index \(2\), and its image already has second index \(p\); the remaining \(c-1\) restrictions are those above. The composition
\(\tau_c\cdots\tau_1\)
is \(\operatorname{tr}_c\) and satisfies
\[
\|\operatorname{tr}_cu\|_{B^{s-\frac{c}{p}}_{p,p}(\mathbb R^k)}
\leq C\|u\|_{H^{s,p}(\mathbb R^n)}.
\]

Starting from \(F^s_{p,q}(\mathbb R^n)\), the first restriction has image \(B^{s-\frac{1}{p}}_{p,p}(\mathbb R^{n-1})=F^{s-\frac{1}{p}}_{p,p}(\mathbb R^{n-1})\), so the preceding calculation applies again. Starting from \(B^s_{p,q}(\mathbb R^n)\), every restriction preserves the index \(q\) and lowers the order by \(\frac{1}{p}\). This gives the other two estimates.

Let \(e_r\) be the codimension-one coretractions, each inserting the variable \(z_r\). In the construction of Theorem~\ref{teo:traza-fraccionaria-euclidea-hiperplano}, a single family of profiles in the normal variable is fixed; thus, the same operators \(e_r\) act simultaneously on all scales in the statement. Choose those operators and define
\[
\operatorname{ex}_c:=e_1e_2\cdots e_c.
\]
The composition is linear and continuous on each of the three scales. Moreover, canceling the factors from the outside inward gives
\[
\operatorname{tr}_c\operatorname{ex}_c
=\tau_c\cdots\tau_1e_1\cdots e_c=I.
\]
Thus, \(\operatorname{tr}_c\) is surjective and \(\operatorname{ex}_c\) is a continuous right inverse.
\end{proof}

For restrictions computed in different charts to represent the same datum, the trace must commute with multiplication by cutoffs and with coordinate changes preserving the submanifold.

\begin{lemma}[Compatibility properties of the Euclidean trace]
\label{lem:compatibilidades-traza-euclidiana}
Let $1<p<\infty$, $1\leq q\leq\infty$, and $s>\frac{c}{p}$. Consider any of the three continuous operators
\[
\begin{aligned}
\operatorname{tr}_c&\colon H^{s,p}(\mathbb R^n)
\longrightarrow B^{s-\frac cp}_{p,p}(\mathbb R^k),\\
\operatorname{tr}_c&\colon F^s_{p,q}(\mathbb R^n)
\longrightarrow B^{s-\frac cp}_{p,p}(\mathbb R^k),\\
\operatorname{tr}_c&\colon B^s_{p,q}(\mathbb R^n)
\longrightarrow B^{s-\frac cp}_{p,q}(\mathbb R^k).
\end{aligned}
\]
The following properties hold in each case.
\begin{enumerate}[label=(\alph*)]
\item Let $a\in C^\infty(\mathbb R^n)$. Suppose multiplication by $a$ is continuous on the source space and multiplication by $a\restriction_{\mathbb R^k}$ is continuous on the target space. In the cases with fine index $q=\infty$, also assume that both multiplications remain continuous when $s$ is replaced by $s-\varepsilon$ and the order of the target space by $s-\varepsilon-\frac{c}{p}$, for every sufficiently small $\varepsilon>0$ with $s-\varepsilon>\frac{c}{p}$. Then
\[
\operatorname{tr}_c(au)
=
(a\restriction_{\mathbb R^k})\operatorname{tr}_cu.
\]
\item Let $\Phi\colon \mathbb R^n\longrightarrow\mathbb R^n$ be a diffeomorphism preserving $\mathbb R^k$. Suppose composition with $\Phi$ and with $\Phi\restriction_{\mathbb R^k}$ is continuous on the source and target spaces, respectively. In the cases with fine index $q=\infty$, also assume the same continuity after replacing $s$ by $s-\varepsilon$ and the target-space order by $s-\varepsilon-\frac{c}{p}$, for every sufficiently small $\varepsilon>0$ with $s-\varepsilon>\frac{c}{p}$. Then
\[
\operatorname{tr}_c(u\circ\Phi)
=
(\operatorname{tr}_cu)\circ
(\Phi\restriction_{\mathbb R^k}).
\]
\item The same identities hold componentwise for functions taking values in a finite-dimensional vector space.
\end{enumerate}
\end{lemma}

\begin{proof}
For $u\in\mathcal S(\mathbb R^n)$, the identities in (a) and (b) are pointwise equalities. First consider $H^{s,p}(\mathbb R^n)$ or one of the spaces $F^s_{p,q}(\mathbb R^n)$ and $B^s_{p,q}(\mathbb R^n)$ with $q<\infty$. In these cases, Corollary~\ref{cor:aproximacion-escalas-BF} gives density of $\mathcal S(\mathbb R^n)$ in the source space. Continuity of the trace, multiplication, and composition permits passage to the limit and proves the identities.

For $F^s_{p,\infty}(\mathbb R^n)$ and $B^s_{p,\infty}(\mathbb R^n)$, fix $\varepsilon\in(0,s-c/p)$ and write $A$ for the scale under consideration. Corollary~\ref{cor:aproximacion-escalas-BF} provides $u_j\in\mathcal S(\mathbb R^n)$ such that
\[
 u_j\longrightarrow u
 \quad\text{in }A^{s-\varepsilon}_{p,\infty}(\mathbb R^n).
\]
For example, first approximate in the space with fine index $p$ and order $s-\varepsilon/2$, and then use its embedding into the displayed space. The loss of order permits this approximation even though Schwartz functions need not be dense in the original norm.

Traces of the same element taken at two admissible orders agree. Indeed, in codimension one they are defined by the same series of block restrictions in Theorem~\ref{teo:traza-fraccionaria-euclidea-hiperplano}; its distributional limit is independent of the order used to estimate it. Composing these traces preserves compatibility in codimension $c$. The estimates place the limit in the target space of order $s-c/p$, even though the approximations converge only at the lower order.

Now apply the pointwise identities of the lemma to $u_j$. The additional continuities in (a) and (b), all at fine index $\infty$ and order $s-\varepsilon$, permit passage to the limit on both sides. Continuity of the trace at that order identifies their limits as distributions on $\mathbb R^k$. By the compatibility just proved, these limits are the traces of order $s$. Both sides belong to the target space of order $s-c/p$, so the equality holds there as well.

For (c), fix a basis of the value space and apply the scalar identities to each component. Since there are only finitely many components, all preceding operations and convergences are preserved.
\end{proof}

\section{Pairs of bounded geometry and uniform Fermi estimates}

In the Fermi charts of this section, write the \(c\) normal variables first, followed by the \(k\) tangential variables. Thus, \(N\) is represented by \(\{0\}^c\times\mathbb R^k\), exactly as in the Euclidean model of the preceding section.

\begin{definition}[Pair of bounded geometry]
\label{def:par-geometria-acotada-traza}
\index{pair of bounded geometry}
Let \((M,\mathbf{g})\) be a Riemannian manifold without boundary of dimension \(n\), and let \(N\subseteq M\) be an embedded submanifold without boundary, closed as a subset of \(M\), of dimension \(k\), with $0\leq k<n$, equipped with the induced metric. Set \(c:=n-k\). We say that \((M,N)\) is a \textit{pair of bounded geometry} if the following conditions hold.
\begin{enumerate}
\item \((M,\mathbf{g})\) has bounded geometry in the sense of Definition~\ref{def: geometria acotada}.
\item \(\operatorname{inj}(N)>0\).
\item There exists \(r_\perp>0\) such that
\[
\exp^\perp\colon
\left\{(x,v)\in T^\perp N\middle| |v|_{\mathbf{g}}<r_\perp\right\}
\longrightarrow M
\]
is a diffeomorphism onto its image.
\item For each \(\ell\in\mathbb N_0\),
\[
\sup_{x\in N}
\bigl|(\nabla^{\mathbf{II}})^\ell\mathbf{II}(x)\bigr|_{\mathbf{g}}<\infty,
\]
where the covariant derivatives and pointwise norm are those of Definition~\ref{def:derivadas-covariantes-segunda-forma-fundamental}.
\end{enumerate}
This is the bounded-geometry condition for the pair \((M,N)\) in \cite[Definition~18]{traces}.
\end{definition}

\begin{remark}
The third condition is global: it prevents normal tubes based at distinct points of \(N\) from meeting at arbitrarily small distance. A uniformly finite local version suffices for continuity of the trace; the global right-inverse construction we will use requires a uniform tubular neighborhood. This distinction corresponds to \cite[Remark~33]{traces}.
\end{remark}

\begin{lemma}[Gauss formula and induced geometry]
\label{lem:subvariedad-geometria-acotada}
Let $(M,\mathbf{g})$ be a Riemannian manifold without boundary, and let $N\subseteq M$ be an embedded submanifold without boundary, closed as a subset of $M$. If \((M,N)\) is a pair of bounded geometry, then \((N,\mathbf{g}\restriction_N)\) is a manifold of bounded geometry. More precisely, for each \(\ell\in\mathbb N_0\), there exists a constant \(C'_\ell\), determined by the bounds on \((\nabla^M)^j\mathbf{R}^M\) and \((\nabla^{\mathbf{II}})^j\mathbf{II}\) with \(j\leq\ell\), such that
\[
\bigl|(\nabla^N)^\ell \mathbf{R}^N\bigr|_{\mathbf{g}\restriction_N}\leq C'_\ell.
\]
\end{lemma}

\begin{proof}
For fields \(\mathbf{X},\mathbf{Y}\) tangent to \(N\), the orthogonal decomposition of the ambient derivative is
\[
\nabla^M_{\mathbf{X}}\mathbf{Y}=\nabla^N_{\mathbf{X}}\mathbf{Y}+\mathbf{II}(\mathbf{X},\mathbf{Y}).
\]
If \(\mathbf{X},\mathbf{Y},\mathbf{Z},\mathbf{W}\) are tangent, substitute this identity into the definition of \(\mathbf{R}^M(\mathbf{X},\mathbf{Y})\mathbf{Z}\), take tangential components, and pair with \(\mathbf{W}\). Terms with two tangential derivatives form \(\mathbf{R}^N\); normal terms are transformed using
\[
\mathbf{g}(\nabla^M_{\mathbf{X}}\boldsymbol{\eta},\mathbf{W})=-\mathbf{g}(\boldsymbol{\eta},\mathbf{II}(\mathbf{X},\mathbf{W})),
\qquad \boldsymbol{\eta}\in\Gamma(T^\perp N),
\]
which follows by differentiating \(\mathbf{g}(\boldsymbol{\eta},\mathbf{W})=0\). This gives the Gauss equation
\begin{equation}
\begin{split}
\mathbf{g}(\mathbf{R}^N(\mathbf{X},\mathbf{Y})\mathbf{Z},\mathbf{W})
={}&\mathbf{g}(\mathbf{R}^M(\mathbf{X},\mathbf{Y})\mathbf{Z},\mathbf{W})\\
&+\mathbf{g}(\mathbf{II}(\mathbf{X},\mathbf{W}),\mathbf{II}(\mathbf{Y},\mathbf{Z}))
-\mathbf{g}(\mathbf{II}(\mathbf{X},\mathbf{Z}),\mathbf{II}(\mathbf{Y},\mathbf{W})).
\end{split}
\label{eq:gauss-traza}
\end{equation}

Apply $(\nabla^N)^\ell$ to \eqref{eq:gauss-traza}, using compatibility with contractions from Proposition~\ref{propiedades de la derivada covariante tensorial} and the Leibniz rule for the operator $*$ from Proposition~\ref{prop: leibniz operador *}. To specify the terms involved, define
\[
 \mathcal I_\ell:=
 \left\{(a,r,b_1,\ldots,b_r)\;\middle|\;
 a,r\in\mathbb N_0,\ (b_1,\ldots,b_r)\in\mathbb N_0^r,
 \ a+\sum_{j=1}^r(b_j+1)\leq\ell\right\}.
\]
The empty sum is zero when $r=0$. We claim that $(\nabla^N)^\ell(\mathbf R^M\restriction_{TN})$ is a finite sum of contractions, with universal coefficients, of terms
\begin{equation}
 \bigl((\nabla^M)^a\mathbf R^M\bigr)\restriction_N
 * (\nabla^{\mathbf{II}})^{b_1}\mathbf{II}
 *\cdots*(\nabla^{\mathbf{II}})^{b_r}\mathbf{II},
 \qquad (a,r,b_1,\ldots,b_r)\in\mathcal I_\ell.
 \label{eq:derivadas-restriccion-curvatura}
\end{equation}
The case $\ell=0$ has a single factor: $a=r=0$. Suppose the expression has been obtained at order $\ell$. When a summand is differentiated, the Leibniz rule makes the connection act on exactly one factor. If it acts on the ambient curvature, the leading term replaces $a$ by $a+1$. Comparing the ambient connection with the tangential and normal connections, the Gauss and Weingarten formulas also produce a factor $\mathbf{II}$: replace $r$ by $r+1$ and add the exponent $b_{r+1}=0$. If the derivative acts on factor number $j\in\{1,\ldots,r\}$, it replaces $b_j$ by $b_j+1$. In all three cases, the weight $\displaystyle a+\displaystyle\sum_{j=1}^r(b_j+1)$ increases by one. Metric contractions produce no additional terms, since the connections preserve the metrics. Thus, all new indices belong to $\mathcal I_{\ell+1}$ and the number of terms still depends only on $n,\dim N,\ell+1$. This proves the expression for every $\ell\in\mathbb N_0$.

The derivatives of the quadratic Gauss term are
\[
(\nabla^{\mathbf{II}})^a\mathbf{II}
*
(\nabla^{\mathbf{II}})^b\mathbf{II},
\qquad a+b=\ell.
\]
The symbol \(*\) denotes a universal contraction with the metric, and the number of summands depends only on \(n,k,\ell\). The hypotheses of Definition~\ref{def:par-geometria-acotada-traza} bound every factor and give \(C'_\ell\). Since \(\operatorname{inj}(N)>0\) is already part of the definition, \(N\) has bounded geometry by Definition~\ref{def: geometria acotada}.
\end{proof}

\begin{lemma}[Uniform estimates for families of differential equations]
\label{lem:edo-uniforme-parametros}
Let $P\subseteq\mathbb R^q$ and $V\subseteq\mathbb R^d$ be open, let $P_0\Subset P$, and let $\Lambda$ be an index set. Consider
\[
\dot y=F_\lambda(t,a,y),
\qquad
y(0)=b_\lambda(a),
\qquad
a\in P_0.
\]
Suppose $F_\lambda$ is continuous in $t$, smooth in $(a,y)$, and all its derivatives with respect to $(a,y)$ are jointly continuous. Suppose $b_\lambda$ is smooth and there exists a compact set $K\Subset V$ with the following property. For all multi-indices $\alpha\in\mathbb N_0^q$ and $\beta\in\mathbb N_0^d$, there exist constants independent of $\lambda$ such that
\[
\sup_{\lambda\in\Lambda}
\sup_{(t,a,y)\in[0,1]\times P_0\times K}
\bigl\|\partial_a^\alpha\partial_y^\beta
F_\lambda(t,a,y)\bigr\|
\leq C_{\alpha,\beta},
\]
\[
\sup_{\lambda\in\Lambda}
\sup_{a\in P_0}
\bigl\|\partial_a^\alpha b_\lambda(a)\bigr\|
\leq C_\alpha'.
\]
Here derivatives are ordinary partial derivatives in finite-dimensional spaces, taken componentwise. If the solutions $y_\lambda(t,a)$ exist for $0\leq t\leq1$ and remain in $K$, then, for each $m\in\mathbb N_0$, there exists $C_m''>0$ such that
\[
\sup_{\lambda\in\Lambda}
\sup_{a\in P_0}
\sup_{t\in[0,1]}
\max_{\substack{\alpha\in\mathbb N_0^q\\|\alpha|\leq m}}
\bigl\|\partial_a^\alpha y_\lambda(t,a)\bigr\|
\leq C_m''.
\]
\end{lemma}

\begin{proof}
Fix $\lambda$ and a parameter $a$ in a neighborhood of $P_0$ where the solution exists. The integral equation is
\[
 y_\lambda(t,a)=b_\lambda(a)
   +\int_0^tF_\lambda(\tau,a,y_\lambda(\tau,a))\,d\tau.
\]
On a compact neighborhood of the trajectory, $\partial_yF_\lambda$ is bounded. Subtracting the equations for two parameters and applying Grönwall first proves continuity of $a\mapsto y_\lambda(\cdot,a)$ in $C([0,1])$. For a coordinate direction, set $Z_h(t):=(y_\lambda(t,a+he_\nu)-y_\lambda(t,a))/h$. The fundamental theorem of calculus in the variables $(a,y)$ gives
\[
 Z_h(t)=\frac{b_\lambda(a+he_\nu)-b_\lambda(a)}h
       +\int_0^t\bigl(A_h(\tau)Z_h(\tau)+f_h(\tau)\bigr)\,d\tau,
\]
where $A_h$ and $f_h$ are averages of $\partial_yF_\lambda$ and $\partial_{a_\nu}F_\lambda$ along the segment joining the two arguments. By joint continuity, they converge uniformly to $\partial_yF_\lambda(\tau,a,y_\lambda)$ and $\partial_{a_\nu}F_\lambda(\tau,a,y_\lambda)$. Applying Grönwall again to the difference between this equation and its limiting equation proves uniform convergence of $Z_h$. Thus, $Z_\nu=\partial_{a_\nu}y_\lambda$ exists and satisfies the variational equation
\[
\dot Z_\nu
=
(\partial_yF_\lambda)(t,a,y_\lambda)Z_\nu
+
\partial_{a_\nu}F_\lambda(t,a,y_\lambda),
\qquad
Z_\nu(0,a)=\partial_{a_\nu}b_\lambda(a).
\]
The hypotheses and Grönwall's Lemma~\ref{lema: gronwall} bound $Z_\nu$ uniformly, simultaneously for $\nu\in\{1,\dots,q\}$.

The same difference-quotient argument, applied to the variational equations already obtained, proves inductively the existence and continuity of all parameter derivatives. For their uniform bounds, suppose derivatives of order less than $m\geq2$ have already been bounded. Fix a multi-index $\gamma\in\mathbb N_0^q$ with $|\gamma|=m$. Apply $\partial_a^\gamma$ to the identity
\[
\partial_t y_\lambda(t,a)
=
F_\lambda(t,a,y_\lambda(t,a)).
\]
The multivariable Faà di Bruno formula of Theorem~\ref{faa di bruno multivariable}, applied componentwise to the function of the variables $(a,y)$ and the mapping $a\mapsto(a,y_\lambda(t,a))$, separates exactly the term in which a derivative with respect to $y$ falls on $\partial_a^\gamma y_\lambda$. This gives
\[
\partial_t(\partial_a^\gamma y_\lambda)
=
(\partial_yF_\lambda)(t,a,y_\lambda)
\partial_a^\gamma y_\lambda
+
Q_{\lambda,\gamma}(t,a),
\]
\[
\partial_a^\gamma y_\lambda(0,a)
=
\partial_a^\gamma b_\lambda(a).
\]
Each summand of $Q_{\lambda,\gamma}$ is a product of a derivative \[
\partial_a^\alpha\partial_y^\beta F_\lambda
(t,a,y_\lambda(t,a))
\] and derivatives $\partial_a^{\gamma_j}y_\lambda$ with $|\gamma_j|<m$. The number of summands depends only on $m,q,d$. The hypotheses and induction therefore provide a constant $K_m$ such that \[
\|Q_{\lambda,\gamma}(t,a)\|\leq K_m
\] for all parameters considered. Another application of Grönwall's Lemma~\ref{lema: gronwall} bounds $\partial_a^\gamma y_\lambda$ uniformly in $[0,1]$. Since there are only finitely many multi-indices of order at most $m$, induction proves the estimate in the statement.
\end{proof}

\begin{remark}
Lemma~\ref{lema: estimacion derivadas flujo variedad} treats a single autonomous field and gives a more precise estimate along each trajectory. The preceding lemma allows time dependence, external parameters, and variable initial data. In particular, dependence on initial data is recovered by taking $F_\lambda$ independent of $a$ and $b_\lambda(a)=a$.
\end{remark}

\begin{proposition}[Uniform Fermi atlas and partition]
\label{prop:atlas-fermi-uniforme}
Let \((M,N)\) be a pair of bounded geometry in the sense of Definition~\ref{def:par-geometria-acotada-traza}; in particular, \(M\) and \(N\) are smooth manifolds without boundary of dimensions \(n\) and \(k\), respectively. There exist \(r_{\mathrm F}>0\), a countable family of points \((p_i)_{i\in I_N}\subseteq N\), Fermi charts
\[
\kappa_i\colon
B_{\mathrm{euc}}^c(0,2r_{\mathrm F})\times B_{\mathrm{euc}}^k(0,2r_{\mathrm F})
\longrightarrow U_i\subseteq M,
\qquad i\in I_N,
\]
geodesic charts \((\kappa_j)_{j\in I_0}\) whose images cover the complement of a smaller tube and satisfy \(U_j\cap N=\varnothing\), and functions \(h_\alpha\in C_c^\infty(U_\alpha)\), \(\alpha\in I:=I_N\sqcup I_0\), with the following properties.
\begin{enumerate}[label=(\alph*)]
\item The cover \((U_\alpha)_{\alpha\in I}\) is uniformly locally finite in the sense of Definition~\ref{def:geometria-acotada-espacio-topologico-cubierta-abierta-uniformemente-localmente}.
\item \(\displaystyle\sum_{\alpha\in I}h_\alpha^2=1\) on \(M\), and, for every multi-index \(\beta\), there exists \(C_\beta\) such that
\[
\sup_{\alpha\in I}
\bigl\|D^\beta(h_\alpha\circ\kappa_\alpha)\bigr\|_{L^\infty(\operatorname{dom}\kappa_\alpha)}
\leq C_\beta.
\]
\item If \(i\in I_N\), then
\begin{equation}
\kappa_i(z,u)
=
\exp^\perp\!\left(
\exp^N_{p_i}(\lambda_i u),
\sum_{a=1}^c z^a\boldsymbol{\nu}_{i,a}(\exp^N_{p_i}(\lambda_i u))
\right),
\label{eq:carta-fermi-uniforme}
\end{equation}
where \(\lambda_i\colon \mathbb R^k\longrightarrow T_{p_i}N\) is an isometry and \((\boldsymbol{\nu}_{i,a})_{a=1}^c\) is an orthonormal frame of the normal bundle obtained by parallel transport for the normal connection along radial geodesics of \(N\). Moreover,
\[
\kappa_i^{-1}(N\cap U_i)
=
\{0\}\times B_{\mathrm{euc}}^k(0,2r_{\mathrm F}).
\]
\item In every chart, including Fermi charts, there exist \(0<c_0\leq C_0<\infty\) such that
\begin{equation}
c_0\|\zeta\|^2
\leq
\sum_{a,b=1}^n g_{\alpha,ab}(x)\zeta^a\zeta^b
\leq
C_0\|\zeta\|^2.
\label{eq:elipticidad-fermi}
\end{equation}
For each multi-index \(\beta\), the functions \(D^\beta g_{\alpha,ab}\) and \(D^\beta \mathbf{g}_\alpha^{ab}\) are bounded by a constant independent of \(\alpha\).
\item Transitions between these charts and the geodesic atlas of Theorem~\ref{thm:trivializacion-geodesica}, and their inverses, have uniformly bounded derivatives of every order on the domains where they are used. Thus, \[
\mathcal Q_{\mathrm F}:=(U_\alpha,\kappa_\alpha^{-1},h_\alpha)_{\alpha\in I}
\] is an admissible quadratic localization system and \((U_\alpha,\kappa_\alpha^{-1},h_\alpha^2)_{\alpha\in I}\) is the associated admissible trivialization.
\end{enumerate}
\end{proposition}

\begin{proof}
Begin by obtaining estimates valid for any center in $N$. By Lemma~\ref{lem:subvariedad-geometria-acotada}, \(N\) has bounded geometry. Let \(r_{0,N}\) and \(r_{0,M}\) be radii for which condition (2) of Theorem~\ref{teo: equivalencias de geometria acotada} holds on \(N\) and \(M\), respectively. Lemma~\ref{lema: volumen bolas geometria acotada} applies up to these radii. Choose
\[
0<r_{\mathrm F}<
\min\left\{
\frac{1}{8}\operatorname{inj}(N),
\frac{r_\perp}{8},
\frac{1}{8}\operatorname{inj}(M),
\frac{r_{0,N}}{5},
\frac{r_{0,M}}{5}
\right\}.
\]
For now, the index $i$ denotes an arbitrary center $p_i\in N$, rather than a previously chosen net. The following estimates are independent of this center. Once all radius reductions have been fixed, we will choose the maximal net and verify the cover multiplicity.

To construct the normal frames, fix an orthonormal basis of \(T^\perp_{p_i}N\) and parallel transport it, for the normal connection, along the radial geodesics of \(N\) starting at \(p_i\). First specify all connections involved. Denote by \(\nabla^M\) and \(\nabla^N\) the Levi--Civita connections of \(M\) and \(N\), by \[
\nabla_{\mathbf{X}}^\perp\boldsymbol{\eta}=(\nabla_{\mathbf{X}}^M\boldsymbol{\eta})^\perp
\] the normal connection, by \(\nabla^{\mathbf{II}}\) the induced product connection on \(T^*N\otimes T^*N\otimes T^\perp N\), and by \(\nabla^{F}\) the connection induced by \(\nabla^M\) on the pullback of \(TM\) under each variational mapping \(F\) below. Our conventions give the Gauss and Weingarten formulas
\begin{equation}
\nabla_{\mathbf{X}}^M \mathbf{Y}=\nabla_{\mathbf{X}}^N \mathbf{Y}+\mathbf{II}(\mathbf{X},\mathbf{Y}),
\qquad
\nabla_{\mathbf{X}}^M\boldsymbol{\eta}
=-\mathbf{A}_{\boldsymbol{\eta}}\mathbf{X}+\nabla_{\mathbf{X}}^\perp\boldsymbol{\eta},
\qquad
\mathbf{g}(\mathbf{A}_{\boldsymbol{\eta}}\mathbf{X},\mathbf{Y})
=\mathbf{g}(\mathbf{II}(\mathbf{X},\mathbf{Y}),\boldsymbol{\eta}).
\label{eq:weingarten-induccion-fermi}
\end{equation}

Taking the normal component in the definition of \(\mathbf{R}^M(\mathbf{X},\mathbf{Y})\boldsymbol{\eta}\) and using \eqref{eq:weingarten-induccion-fermi} gives the Ricci equation
\begin{equation}
\mathbf{g}(\mathbf{R}^\perp(\mathbf{X},\mathbf{Y})\boldsymbol{\eta},\boldsymbol{\zeta})
=
\mathbf{g}(\mathbf{R}^M(\mathbf{X},\mathbf{Y})\boldsymbol{\eta},\boldsymbol{\zeta})
+\mathbf{g}([\mathbf{A}_{\boldsymbol{\eta}},\mathbf{A}_{\boldsymbol{\zeta}}]\mathbf{X},\mathbf{Y}).
\label{eq:ricci-conexion-normal}
\end{equation}
Indeed, the normal terms of \(\nabla_{\mathbf{X}}^M(-\mathbf{A}_{\boldsymbol{\eta}}\mathbf{Y})\) and \(\nabla_{\mathbf{Y}}^M(-\mathbf{A}_{\boldsymbol{\eta}}\mathbf{X})\) are \(-\mathbf{II}(\mathbf{X},\mathbf{A}_{\boldsymbol{\eta}}\mathbf{Y})\) and \(-\mathbf{II}(\mathbf{Y},\mathbf{A}_{\boldsymbol{\eta}}\mathbf{X})\); pairing them with \(\boldsymbol{\zeta}\) gives the commutator of \eqref{eq:ricci-conexion-normal}. Derivatives of this identity are controlled by the following recurrence. At order zero, take exactly the two summands on its right-hand side. To pass from order $q$ to $q+1$, apply the product connection to each contraction and use the exact rule
\[
 \nabla_{\mathbf{X}}\operatorname{contr}(\mathbf{T}_1\otimes\cdots\otimes \mathbf{T}_\ell)
 =\sum_{j=1}^{\ell}
 \operatorname{contr}(\mathbf{T}_1\otimes\cdots\otimes
 \nabla_{\mathbf{X}}\mathbf{T}_j\otimes\cdots\otimes \mathbf{T}_\ell).
\]
If the derivative reaches the restriction of $(\nabla^M)^a\mathbf{R}^M$, the Gauss and Weingarten formulas produce $(\nabla^M)^{a+1}\mathbf{R}^M$ and contractions of $(\nabla^M)^a\mathbf{R}^M$ with $\mathbf{II}$. If it reaches a shape operator, the identity
\[
 \mathbf{g}\bigl((\nabla_{\mathbf{X}}\mathbf{A})_{\boldsymbol{\eta}}\mathbf{Y},\mathbf{Z}\bigr)
 =\mathbf{g}\bigl((\nabla_{\mathbf{X}}^{\mathbf{II}}\mathbf{II})(\mathbf{Y},\mathbf{Z}),\boldsymbol{\eta}\bigr)
\]
---understood with the product connections on all arguments---replaces it by a derivative of $\mathbf{II}$. This rule determines a finite family at each order and proves by induction that $(\nabla^\perp)^q\mathbf{R}^\perp$ is a universal sum of contractions of $(\nabla^M)^a\mathbf{R}^M$ and $(\nabla^{\mathbf{II}})^b\mathbf{II}$, with $a,b\leq q$, together with undifferentiated copies of $\mathbf{II}$ introduced by Gauss--Weingarten.

If \(\omega^\perp_\mu\) and \(F^\perp_{\nu\mu}\) are the matrices of the normal connection and its curvature, respectively, in a radial geodesic chart of \(N\), radial gauge satisfies \(\displaystyle \sum_{\mu=1}^k y^\mu\omega^\perp_\mu(y)=0\) and the exact identity
\begin{equation}
\omega^\perp_\mu(y)
=
\int_0^1 s\sum_{\nu=1}^k y^\nu F^\perp_{\nu\mu}(sy)\,ds.
\label{eq:calibre-radial-normal-fermi}
\end{equation}
This follows by contracting the structure equation with the radial field and differentiating \(s\omega^\perp_\mu(sy)\). For a multi-index $\alpha$, the Leibniz rule and linearity of $y^\nu$ give the full formula
\begin{align*}
 \partial^\alpha\omega^\perp_\mu(y)
 =\int_0^1{}&s^{|\alpha|+1}\sum_{\nu=1}^k y^\nu
 (\partial^\alpha F^\perp_{\nu\mu})(sy)\,ds\\
 &+\sum_{\substack{j\in\{1,\ldots,k\}\\\alpha_j>0}}\alpha_j
 \int_0^1s^{|\alpha|}
 (\partial^{\alpha-e_j}F^\perp_{j\mu})(sy)\,ds.
\end{align*}
Control is established by simultaneous induction. At order zero, the tensor bound on $\mathbf{R}^\perp$ bounds its components $F^\perp_{\nu\mu}$ in the orthonormal frame, and the radial formula bounds $\omega^\perp$. Suppose derivatives of $\omega^\perp$ of order less than $q$ have been controlled. The component identity
\[
 (\nabla^\perp_\lambda F^\perp)_{\nu\mu}
 =\partial_\lambda F^\perp_{\nu\mu}
  +[\omega^\perp_\lambda,F^\perp_{\nu\mu}]
  -\sum_{\kappa=1}^k\left(
   \Gamma^\kappa_{\lambda\nu}F^\perp_{\kappa\mu}
  +\Gamma^\kappa_{\lambda\mu}F^\perp_{\nu\kappa}\right)
\]
and the higher-derivative expression in Lemma~\ref{lem:local-expression-higher-order} express every ordinary derivative of $F^\perp$ of order $q$ through $(\nabla^\perp)^aF^\perp$, $a\leq q$, derivatives of $\omega^\perp$ of order less than $q$, and already controlled derivatives of the Christoffel symbols of $N$. The preceding integral formula then bounds derivatives of $\omega^\perp$ of order $q$ and closes the induction.

Frame bounds are interpreted through the connection, not through components in an arbitrary trivialization. Indeed,
\[
 \nabla^\perp_{\boldsymbol\partial_\mu}\boldsymbol\nu_{i,a}
 =\sum_{b=1}^c(\omega^\perp_\mu)^b{}_a\boldsymbol\nu_{i,b}.
\]
Since the frame is orthonormal, this formula bounds its first covariant derivatives; differentiating it, the Leibniz rule and the bounds on $\partial^\alpha\omega^\perp$ bound derivatives of every order. Gauss--Weingarten converts these normal controls into controls for the ambient connection, adding only factors of $\mathbf{II}$ and its derivatives. In uniform ambient charts, the local expression for the connection finally converts these controls into component bounds. These are the initial data used in the Jacobi equation.

Now consider the normal exponential map as the mapping
\[
\mathcal F\colon
\mathcal V_\perp
:=
\{(x,\xi)\in T^\perp N\mid |\xi|_{\mathbf{g}}<r_\perp\}
\longrightarrow M,
\qquad
\mathcal F(x,\xi)=\exp_x^M(\xi).
\]
By the tubular hypothesis, its restriction to \(\mathcal V_\perp\) is a diffeomorphism onto its image. Set \[
x_i(u):=\exp^N_{p_i}(\lambda_i u),
\qquad
\boldsymbol{\xi}_i(z,u):=\sum_{a=1}^c z^a\boldsymbol{\nu}_{i,a}(x_i(u)),
\] and, for \(0\leq t\leq1\),
\begin{equation}
F_i(z,u,t):=\mathcal F(x_i(u),t\boldsymbol{\xi}_i(z,u)).
\label{eq:familia-geodesica-fermi}
\end{equation}
Then \(\kappa_i(z,u)=F_i(z,u,1)\).

Write \[
D:=\nabla^{F_i}_{\boldsymbol{\partial}_t},\qquad
\nabla_{\mathbf{n}_a}:=\nabla^{F_i}_{\boldsymbol{\partial}_{z^a}},\qquad
\nabla_{\mathbf{t}_\mu}:=\nabla^{F_i}_{\boldsymbol{\partial}_{u^\mu}},
\] and use the ordered alphabet
\[
\mathscr A
=
\{\mathbf{n}_1,\ldots,\mathbf{n}_c,
  \mathbf{t}_1,\ldots,\mathbf{t}_k\}.
\]
Let
\[
\mathbf{T}:=dF_i(\boldsymbol{\partial}_t),\qquad
\mathbf{J}_{\mathfrak a}:=dF_i(\boldsymbol{\partial}_{\mathfrak a}),
\qquad \mathfrak a\in\mathscr A.
\]
Since the domain's coordinate fields commute and \(\nabla^M\) is torsion-free,
\begin{equation}
D\mathbf{T}=0,\qquad
\nabla_{\mathfrak a}\mathbf{T}=D\mathbf{J}_{\mathfrak a},\qquad
\nabla_{\mathfrak a}D\mathbf{Z}-D\nabla_{\mathfrak a}\mathbf{Z}
=
\mathbf{R}^M(\mathbf{J}_{\mathfrak a},\mathbf{T})\mathbf{Z}.
\label{eq:conmutadores-fermi}
\end{equation}
Each \(\mathbf{J}_{\mathfrak a}\) satisfies the Jacobi equation
\[
D^2\mathbf{J}_{\mathfrak a}
+\mathbf{R}^M(\mathbf{J}_{\mathfrak a},\mathbf{T})\mathbf{T}=0.
\]

Fix words with their exact order retained. For a nonempty word \(I=(\mathfrak a_1,\ldots,\mathfrak a_m)\), define
\begin{equation}
\mathbf{Y}_I
:=
\nabla_{\mathfrak a_1}\cdots
\nabla_{\mathfrak a_{m-1}}\mathbf{J}_{\mathfrak a_m}.
\label{eq:palabras-ordenadas-fermi}
\end{equation}
The letters of \(I\) are not permuted. Introduce the operator
\[
\mathcal L \mathbf{Z}:=D^2\mathbf{Z}+\mathbf{R}^M(\mathbf{Z},\mathbf{T})\mathbf{T}.
\]
From \eqref{eq:conmutadores-fermi}, \(D\mathbf{T}=0\), and the covariant Leibniz rule, we obtain
\begin{equation}
\mathcal L(\nabla_{\mathfrak a}\mathbf{Z})
=
\nabla_{\mathfrak a}(\mathcal L\mathbf{Z})
+\mathscr C_{\mathfrak a}(\mathbf{Z}),
\label{eq:conmutacion-jacobi-fermi}
\end{equation}
where
\begin{align}
\mathscr C_{\mathfrak a}(\mathbf{Z})
={}&
-2\mathbf{R}^M(\mathbf{J}_{\mathfrak a},\mathbf{T})D\mathbf{Z}
-(\nabla_{\mathbf{T}}^M \mathbf{R}^M)(\mathbf{J}_{\mathfrak a},\mathbf{T})\mathbf{Z}
-\mathbf{R}^M(D\mathbf{J}_{\mathfrak a},\mathbf{T})\mathbf{Z}
\nonumber\\
&-(\nabla_{\mathbf{J}_{\mathfrak a}}^M \mathbf{R}^M)(\mathbf{Z},\mathbf{T})\mathbf{T}
-\mathbf{R}^M(\mathbf{Z},D\mathbf{J}_{\mathfrak a})\mathbf{T}
-\mathbf{R}^M(\mathbf{Z},\mathbf{T})D\mathbf{J}_{\mathfrak a}.
\label{eq:fuente-conmutacion-fermi}
\end{align}
In particular, the sources \(\mathbf{Q}_I\) are defined by the finite recurrence, compatible with the tensor spaces,
\begin{equation}
\mathbf{Q}_{(\mathfrak a)}:=0,
\qquad
\mathbf{Q}_{(\mathfrak a,I)}
:=
\nabla_{\mathfrak a}\mathbf{Q}_I+\mathscr C_{\mathfrak a}(\mathbf{Y}_I),
\qquad
\mathcal L \mathbf{Y}_I=\mathbf{Q}_I.
\label{eq:recurrencia-fuentes-fermi}
\end{equation}
The recurrence is closed and indexes every term. Consider formal expressions obtained by contracting, using \(\mathbf{g}\), factors \((\nabla^M)^q\mathbf{R}^M\), \(\mathbf{T}\), \(\mathbf{Y}_J\), and \(D\mathbf{Y}_J\), with a single free vector index. Define the formal derivation \(\mathfrak d_{\mathfrak a}^{F_i}\) by the Leibniz rule on products and contractions and by
\[
\mathfrak d_{\mathfrak a}^{F_i} \mathbf{T}=D\mathbf{J}_{\mathfrak a},
\qquad
\mathfrak d_{\mathfrak a}^{F_i} \mathbf{Y}_J=\mathbf{Y}_{(\mathfrak a,J)},
\qquad
\mathfrak d_{\mathfrak a}^{F_i}(D\mathbf{Y}_J)
=D\mathbf{Y}_{(\mathfrak a,J)}
+\mathbf{R}^M(\mathbf{J}_{\mathfrak a},\mathbf{T})\mathbf{Y}_J,
\]
\[
\mathfrak d_{\mathfrak a}^{F_i}\bigl((\nabla^M)^q\mathbf{R}^M\bigr)
=
\nabla^M_{\mathbf{J}_{\mathfrak a}}\bigl((\nabla^M)^q\mathbf{R}^M\bigr)
=
(\nabla^M)^{q+1}\mathbf{R}^M(\mathbf{J}_{\mathfrak a},\,\cdot\,).
\]
In the last formula, $\cdot$ represents the preceding arguments in their original order, and the new argument is placed first; remaining arguments are differentiated separately by the Leibniz rule. Metric compatibility prevents additional summands when contractions are differentiated. Let \(\mathfrak C^{F_i}_{\mathfrak a,I}\) be the multiset of the six summands, with their coefficients, of \(\mathscr C_{\mathfrak a}(\mathbf{Y}_I)\) in \eqref{eq:fuente-conmutacion-fermi}. Define finite multisets by
\[
\mathfrak A^{F_i}_{(\mathfrak a)}:=\varnothing,
\qquad
\mathfrak A^{F_i}_{(\mathfrak a,I)}
:=
\mathfrak d_{\mathfrak a}^{F_i}\mathfrak A^{F_i}_I
\sqcup\mathfrak C^{F_i}_{\mathfrak a,I}.
\]
Here differentiation of a multiset includes all summands produced by the Leibniz rule and preserves coefficients and multiplicities. Induction on \eqref{eq:recurrencia-fuentes-fermi} gives the exact identity
\begin{equation}
\mathbf{Q}_I=\sum_{\mathfrak b\in\mathfrak A^{F_i}_I}\mathfrak b.
\label{eq:indexacion-exacta-fuentes-fermi}
\end{equation}

Assign weight zero to \(\mathbf{T}\) and weight \(|J|\) to every occurrence of \(\mathbf{Y}_J\) or \(D\mathbf{Y}_J\). For \(|I|=m\geq2\), every summand of \eqref{eq:indexacion-exacta-fuentes-fermi} satisfies
\begin{equation}
\sum_{\text{occurrences of }\mathbf{Y}_J\text{ or }D\mathbf{Y}_J}|J|=m.
\label{eq:peso-exacto-fuentes-fermi}
\end{equation}
Moreover, it contains one or more factors \((\nabla^M)^q\mathbf{R}^M\), with \(0\leq q\leq m-1\); each field carries at most one derivative \(D\), and every field in the source has length at most \(m-1\). Indeed, the six summands of \(\mathscr C_{\mathfrak a}(\mathbf{Y}_I)\) have weight \(1+|I|\), and every rule of \(\mathfrak d_{\mathfrak a}^{F_i}\) increases the weight by one, never produces \(D^2\mathbf{Y}_J\), and increases the order of a curvature derivative by at most one. In particular, the equation of length \(m\) contains no field of length \(m\) outside the unknown term \(\mathbf{Y}_I\).

The initial data are also determined by an exact recurrence. Let
\[
\mathbf{X}_\mu:=dx_i(\boldsymbol{\partial}_{u^\mu}),\qquad
\mathbf{H}_I:=\mathbf{Y}_I|_{t=0},\qquad \mathbf{K}_I:=D\mathbf{Y}_I|_{t=0}.
\]
For words of length one,
\begin{align}
&\mathbf{H}_{(\mathbf{n}_a)}=0,
&&\mathbf{K}_{(\mathbf{n}_a)}=\boldsymbol{\nu}_{i,a}(x_i(u)),
\label{eq:datos-normales-fermi}\\
&\mathbf{H}_{(\mathbf{t}_\mu)}=\mathbf{X}_\mu,
&&\mathbf{K}_{(\mathbf{t}_\mu)}
=\nabla^M_{\mathbf{X}_\mu}\boldsymbol{\xi}_i
=-\mathbf{A}_{\boldsymbol{\xi}_i}\mathbf{X}_\mu+\nabla^\perp_{\mathbf{X}_\mu}\boldsymbol{\xi}_i.
\label{eq:datos-tangenciales-fermi}
\end{align}
The last equality fixes the sign and displays separately the second fundamental form and the normal connection. In our coordinates,
\[
\mathbf{K}_{(\mathbf{t}_\mu)}
=
-\sum_{a=1}^c z^a\mathbf{A}_{\boldsymbol{\nu}_{i,a}}\mathbf{X}_\mu
+\sum_{a=1}^c z^a\nabla^\perp_{\mathbf{X}_\mu}\boldsymbol{\nu}_{i,a}.
\]
The second summand may be nonzero: radial gauge makes it vanish only in the radial direction of the chart of \(N\).

Denote by \(\nabla^0_{\mathfrak a}\) the induced connection on \(t=0\): for \(\mathfrak a=\mathbf{t}_\mu\), it is the derivative \(\nabla^M_{\mathbf{X}_\mu}\), whereas for \(\mathfrak a=\mathbf{n}_a\), it is ordinary differentiation with respect to \(z^a\) in the fixed space \(T_{x_i(u)}M\). For every word \(I\),
\begin{equation}
\mathbf{H}_{(\mathfrak a,I)}
=
\nabla^0_{\mathfrak a}\mathbf{H}_I,
\qquad
\mathbf{K}_{(\mathfrak a,I)}
=
\nabla^0_{\mathfrak a}\mathbf{K}_I
-\mathbf{R}^M(\mathbf{H}_{(\mathfrak a)},\boldsymbol{\xi}_i)\mathbf{H}_I.
\label{eq:recurrencia-datos-fermi}
\end{equation}
The second identity is precisely \(D\nabla_{\mathfrak a}
=\nabla_{\mathfrak a}D-\mathbf{R}^M(\mathbf{J}_{\mathfrak a},\mathbf{T})\) evaluated at \(t=0\).

Formulas \eqref{eq:datos-normales-fermi}--\eqref{eq:recurrencia-datos-fermi}, expanded recursively using \eqref{eq:weingarten-induccion-fermi}, express all \(\mathbf{H}_I,\mathbf{K}_I\) as finite universal sums of contractions of
\[
(\nabla^M)^q\mathbf{R}^M,
\qquad
(\nabla^{\mathbf{II}})^q\mathbf{II},
\qquad
(\nabla^\perp)^q\boldsymbol{\nu}_{i,a},
\qquad
(\nabla^N)^q\mathbf{X}_\mu,
\]
together with powers of \(z\). Thus, ambient curvature appears in the sources \(\mathbf{Q}_I\), whereas \(\mathbf{II}\) and the normal connection appear precisely in the initial data. The bounded-geometry bounds for \(M\) and \(N\), the hypotheses on \(\mathbf{II}\), and \eqref{eq:calibre-radial-normal-fermi} bound all these data uniformly.

We now prove the bounds by induction on word length. Since \(\mathbf{T}\) is parallel and \(|\mathbf{T}|_{\mathbf{g}}=|\boldsymbol{\xi}_i|_{\mathbf{g}}\leq2r_{\mathrm F}\), in a parallel frame along \(t\mapsto F_i(z,u,t)\), equation \(\mathcal L\mathbf{Y}_I=\mathbf{Q}_I\) is a first-order linear system for \((\mathbf{Y}_I,D\mathbf{Y}_I)\) with uniformly bounded matrix. For \(m=1\), the source is zero and the data \eqref{eq:datos-normales-fermi}--\eqref{eq:datos-tangenciales-fermi} are uniformly bounded. If words of length less than \(m\) are controlled, \eqref{eq:recurrencia-fuentes-fermi} bounds \(\mathbf{Q}_I\) for \(|I|=m\), and \eqref{eq:recurrencia-datos-fermi} bounds \(\mathbf{H}_I,\mathbf{K}_I\). Grönwall's Lemma~\ref{lema: gronwall} gives, with a constant independent of \(i,z,u\),
\[
\sup_{t\in[0,1]}
\bigl(|\mathbf{Y}_I(t)|_{\mathbf{g}}+|D\mathbf{Y}_I(t)|_{\mathbf{g}}\bigr)
\leq
C_m\left(
|\mathbf{H}_I|_{\mathbf{g}}+|\mathbf{K}_I|_{\mathbf{g}}+\int_0^1|\mathbf{Q}_I(s)|_{\mathbf{g}}\,ds
\right).
\]
More precisely, $C_m$ depends only on the dimension, $r_{\mathrm F}$, the bounds on $(\nabla^M)^q\mathbf{R}^M$ for $q\leq m-1$, the bounds on $(\nabla^{\mathbf{II}})^q\mathbf{II}$ and the normal connection appearing in the data of order $m$, and the constants already constructed for words of length less than $m$; it does not depend on the chart index or $(z,u)$. This explicitly completes the induction for every tangential and normal word.

We pass from covariant to ordinary derivatives without identifying the two. Let \(\mathbf{P}_1,\ldots,\mathbf{P}_n\) be an orthonormal frame on \(t=0\), formed from a controlled tangential frame of \(N\) and \(\boldsymbol{\nu}_{i,1},\ldots,\boldsymbol{\nu}_{i,c}\), and extended by \(D\mathbf{P}_r=0\). For an ordered word $I=(\mathfrak a_1,\ldots,\mathfrak a_m)$, set
\[
\partial_I:=\partial_{\mathfrak a_1}\cdots\partial_{\mathfrak a_m},
\qquad
\nabla_I:=\nabla_{\mathfrak a_1}\cdots\nabla_{\mathfrak a_m},
\qquad
\partial_\varnothing=\nabla_\varnothing=\operatorname{id}.
\]
The first letter is the outer operator; this word convention does not change the order of the tensor factors of $\nabla^m$. Also set
\[
\mathbf{P}_{r,I}:=
\nabla_{\mathfrak a_1}\cdots\nabla_{\mathfrak a_m}\mathbf{P}_r.
\]
These words satisfy the triangular recurrence
\[
D\mathbf{P}_{r,\varnothing}=0,
\qquad
D\mathbf{P}_{r,(\mathfrak a,I)}
=
\nabla_{\mathfrak a}(D\mathbf{P}_{r,I})
-\mathbf{R}^M(\mathbf{J}_{\mathfrak a},\mathbf{T})\mathbf{P}_{r,I}.
\]
Their data at \(t=0\) are controlled by the frames of \(TN\) and \(T^\perp N\); the preceding induction and Grönwall's Lemma~\ref{lema: gronwall} bound all \(\mathbf{P}_{r,I}\). Finally, if \(I_S\) denotes the ordered subword of \(I\) whose positions belong to \(S\), the iterated Leibniz rule gives the exact identity
\begin{equation}
\partial_I\langle \mathbf{Z},\mathbf{P}_r\rangle
=
\sum_{S\subseteq\{1,\ldots,|I|\}}
\left\langle
\nabla_{I_S}\mathbf{Z},\mathbf{P}_{r,I_{S^c}}
\right\rangle .
\label{eq:paso-triangular-fermi}
\end{equation}
The term with \(S=\{1,\ldots,|I|\}\) contains the highest-order covariant derivative; the others contain lower orders. This triangular relation proves bounds for all ordinary derivatives of the variational fields and, using the geodesic coordinates of $M$ already controlled by Theorem~\ref{teo: equivalencias de geometria acotada}, for all ordinary derivatives of \(F_i\) and \(\kappa_i=F_i(\,\cdot\,,1)\).

On \(z=0\),
\[
d\kappa_i(\boldsymbol{\partial}_{z^a})=\boldsymbol{\nu}_{i,a},
\qquad
d\kappa_i(\boldsymbol{\partial}_{u^\mu})=\mathbf{X}_\mu.
\]
The two subspaces are orthogonal; the normal block is orthonormal, and the tangential block is the differential of a controlled geodesic chart of \(N\). Thus, there exists \(\sigma_0>0\), independent of \(i,u\), such that the smallest singular value of \(D\kappa_i(0,u)\) is at least \(\sigma_0\). The uniform bound for \(D^2\kappa_i\) gives
\[
\|D\kappa_i(z,u)-D\kappa_i(0,u)\|
\leq C|z|.
\]
Decreasing \(r_{\mathrm F}\) once, impose \(2Cr_{\mathrm F}\leq\frac{\sigma_0}{2}\); then \[
\frac{\sigma_0}{2}|v|
\leq |D\kappa_i(z,u)v|
\leq C|v|
\] on all domains.

Bounds for the inverses follow by an explicit induction. If \(\rho\) is a coordinate representation of \(\kappa_i\) and \(\lambda=\rho^{-1}\), order one is already controlled by
\[
 D\lambda(y)=D\rho(\lambda(y))^{-1}.
\]
Suppose derivatives of \(\lambda\) of orders less than \(m\) are uniformly bounded, and fix a multi-index \(\alpha\) with \(|\alpha|=m\geq2\). Apply the Faà di Bruno formula of Theorem~\ref{faa di bruno multivariable} componentwise to the identity \(\rho\circ\lambda=\operatorname{id}\). For each component \(a\), we obtain
\[
 0=D^\alpha(\rho^a\circ\lambda)(y)
 =\sum_{b=1}^n
   (\partial_b\rho^a)(\lambda(y))D^\alpha\lambda^b(y)
   +R^a_\alpha(y),
\]
where \(R^a_\alpha\) is a finite sum of products of the form
\[
 c\,(D^\beta\rho^a)(\lambda(y))
 \prod_{\ell=1}^{|\beta|}D^{\gamma_\ell}
       \lambda^{j_\ell}(y),
 \qquad
 2\leq|\beta|\leq m,
 \quad
 \sum_{\ell=1}^{|\beta|}\gamma_\ell=\alpha,
 \quad |\gamma_\ell|>0.
\]
Since \(|\beta|\geq2\), each \(|\gamma_\ell|\leq m-1\); by the induction hypothesis, all factors containing derivatives of \(\lambda\) are uniformly bounded. The derivative bounds already obtained for \(\rho\) then bound \(R_\alpha=(R^a_\alpha)_{a=1}^n\). In matrix form, the preceding equality is
\[
 D\rho(\lambda(y))D^\alpha\lambda(y)=-R_\alpha(y),
\]
and hence
\[
 D^\alpha\lambda(y)
 =-D\rho(\lambda(y))^{-1}R_\alpha(y).
\]
The uniform bound for \(D\rho^{-1}\) gives the bound of order \(m\) and completes the induction.

On the overlap of two Fermi charts, uniqueness of the normal-exponential representation gives a more precise formula than a mere composition. If
\[
\tau_{ji}:=x_j^{-1}\circ x_i,
\qquad
\boldsymbol{\nu}_{i,a}(x_i(u))
=
\sum_{b=1}^c O_{ji}(u)^b{}_a\,
\boldsymbol{\nu}_{j,b}(x_j(\tau_{ji}(u))),
\]
then
\(O_{ji}(u)\in O(c)\)
and
\begin{equation}
\kappa_j^{-1}\circ\kappa_i(z,u)
=
\bigl(O_{ji}(u)z,\tau_{ji}(u)\bigr).
\label{eq:transicion-fermi-exacta}
\end{equation}
In particular, the transition preserves the zero section and normal fibers and is linear in the normal variable. Moreover,
\[
O_{ji}(u)^b{}_a
=
\mathbf{g}\bigl(\boldsymbol{\nu}_{j,b},\boldsymbol{\nu}_{i,a}\bigr)_{x_i(u)},
\]
so \eqref{eq:paso-triangular-fermi} and the bounds for both frames control all derivatives of \(O_{ji}\); its inverse is its transpose. Transitions between Fermi charts and geodesic charts of $M$, and their inverses, are controlled by the Faà di Bruno formula of Theorem~\ref{faa di bruno multivariable} and the preceding induction.

The Fermi metric coefficients are
\[
g_{AB}=\mathbf{g}(\mathbf{J}_A,\mathbf{J}_B)\big|_{t=1}.
\]
The Leibniz rule and \eqref{eq:paso-triangular-fermi} bound all their derivatives. The two-sided bound for \(D\kappa_i\) proves \eqref{eq:elipticidad-fermi}. Differentiating \(G^{-1}G=I\) gives, for \(\alpha\neq0\),
\[
(\partial^\alpha G^{-1})G
=
-\sum_{\substack{\beta\in\mathbb N_0^n\\\beta\leq\alpha,\ \beta\neq\alpha}}
\binom{\alpha}{\beta}
(\partial^\beta G^{-1})
(\partial^{\alpha-\beta}G),
\]
which closes the inverse-metric bounds by induction. Likewise,
\[
0<c\leq |\det D\kappa_i|\leq C,
\qquad
0<c\leq\sqrt{\det(\mathbf{g})}\leq C,
\]
and all derivatives of both Jacobians are bounded: the determinant is polynomial in the entries, and the square root is differentiated on a compact interval contained in \((0,\infty)\).

Now fix $r:=r_{\mathrm F}$ after all preceding reductions; take it also small enough that $20r<r_{0,M}$ and $10r<r_{0,N}$. Choose a maximal family of centers $p_i\in N$ whose intrinsic balls $B_N(p_i,r/2)$ are disjoint. The balls $B_N(p_i,r)$ cover $N$: otherwise, another disjoint ball could be added. The family is countable because $N$ is second countable. Its multiplicity, even for balls of radius $4r$, is bounded by comparing volumes of the disjoint balls with the volume of a ball of radius $5r$.

Denote by $U_r(N)$ the image of normal vectors of length less than $r$. Also choose a maximal $\rho$-separated net in $M\setminus U_r(N)$, where $0<\rho<r/8$. Its balls of radius $\rho$ cover this complement, and we use balls of radius $2\rho$ as charts. These do not meet $N$. Indeed, since $N$ is closed and $M$ complete, every point at distance less than $r$ from $N$ has a minimizing segment to $N$ orthogonal to $N$ and belongs to $U_r(N)$. The Fermi charts cover $U_{2r}(N)$; consequently, the two families cover $M$ and leave an overlap band.

To compare the two families, count overlaps through volumes of tubular sets over $N$ in the ambient manifold. Use the disjoint sets
\[
 D_i:=\exp^\perp\left(
 \{(x,v)\mid x\in B_N(p_i,r/2),\ |v|<r/2\}\right).
\]
Injectivity of the tube makes them disjoint. The Fermi Jacobian bounds and volume bounds on $N$ give, with a uniform constant,
\[
 \lambda_{\mathbf g}(D_i)
 \geq c\,r^c\lambda_{\mathbf g|_N}(B_N(p_i,r/2))
 \geq c' r^n.
\]
If a Fermi chart meets an ambient ball $B_M(x,2\rho)$, its center is at distance less than $4r+2\rho$ from $x$ and $D_i\subset B_M(x,6r)$. The volume of this last ball does not exceed $Cr^n$. By disjointness, the number of such Fermi charts is at most $C/c'$. In the reverse direction, the disjoint balls of radius $\rho/2$ based at interior centers whose charts meet a Fermi chart lie in an ambient ball of radius $7r$; each has volume at least $c\rho^n$. This bounds their number by $C(r/\rho)^n$, a fixed constant. The same reasoning compares the Fermi family with any uniform geodesic net at the fixed scale. Between Fermi charts, the normal footpoint of an overlap point lies in both balls of $N$, so their centers are at distance less than $4r$ in $N$. Between interior charts, use geodesic packing. This proves uniform local finiteness and both cross-overlap bounds.

Finally, choose model functions \[
\eta_\perp\in C_c^\infty(B^c(0,2r_{\mathrm F})),
\qquad
\eta_N\in C_c^\infty(B^k(0,2r_{\mathrm F})),
\] equal to one on the balls of radius \(r_{\mathrm F}\), and use \(\eta_\perp(z)\eta_N(u)\) in Fermi charts; in charts away from the submanifold, use a fixed geodesic cutoff. Let \(\widetilde{\eta}_\alpha\) be their transports to \(M\). The covering property and uniform multiplicity give
\[
1\leq
S:=\sum_{\alpha\in I}\widetilde{\eta}_\alpha^2
\leq C.
\]
Only uniformly finitely many summands contribute in any chart. The transition bounds and the Faà di Bruno formula of Theorem~\ref{faa di bruno multivariable} control all derivatives of \(S\) and \(S^{-\frac{1}{2}}\). Thus,
\[
h_\alpha:=\widetilde{\eta}_\alpha S^{-\frac{1}{2}}
\]
satisfies \(\displaystyle \sum_{\alpha\in I} h_\alpha^2=1\) and all estimates in (b). To construct auxiliary cutoffs without losing uniformity, choose once and for all functions \[
 \widehat{\eta}_\perp\in C_c^\infty(B^c(0,2r_{\mathrm F})),
 \qquad
 \widehat{\eta}_N\in C_c^\infty(B^k(0,2r_{\mathrm F})),
\] with values in \([0,1]\), equal to one on neighborhoods of \(\operatorname{supp}\eta_\perp\) and \(\operatorname{supp}\eta_N\), respectively. In a Fermi chart, transport the product \(\widehat{\eta}_\perp(z)\widehat{\eta}_N(u)\); in the other charts, likewise use a single geodesic profile equal to one near the support of the profile used for \(\widetilde{\eta}_\alpha\). The resulting functions \(\chi_\alpha\) belong to \(C_c^\infty(U_\alpha)\), satisfy \(0\leq\chi_\alpha\leq1\), and equal one on a neighborhood of \(\operatorname{supp}h_\alpha\), since \(\operatorname{supp}h_\alpha\subseteq
\operatorname{supp}\widetilde{\eta}_\alpha\). Because only these fixed profiles were used, their coordinate derivatives have bounds independent of \(\alpha\); these are exactly the estimates in \ref{item:B3}.

We now have adapted charts, frames, and cutoffs with uniform constants. In the next section, we can apply the Euclidean trace and sum its estimates over this cover.
\end{proof}
\section{Trace on an interior submanifold}

Let \((U_\alpha,\kappa_\alpha^{-1},h_\alpha)_{\alpha\in I}\) be the quadratic localization system of Proposition~\ref{prop:atlas-fermi-uniforme}. For \(i\in I_N\), write
\[
\kappa_i^N(u):=\kappa_i(0,u),
\qquad
U_i^N:=U_i\cap N,
\qquad
h_i^N:=h_i\restriction_N.
\]
The charts \(\kappa_i^N\) are geodesic charts of \(N\), charts indexed by \(I_0\) do not meet \(N\), and
\begin{equation}
\displaystyle\sum_{i\in I_N}(h_i^N)^2=1
\qquad\text{on }N.
\label{eq:particion-cuadrada-N}
\end{equation}
By Theorems~\ref{teo:geometria-acotada-variedad-riemanniana-geometria-acotada-trivializacion-admisible} and \ref{teo:geometria-acotada-variedad-riemanniana-trivializacion-admisible}, together with Lemma~\ref{lem:equivalencia-localizacion-cuadratica}, for each $s,t\in\mathbb R$ and $1<p<\infty$ there exist constants $c_M,C_M,c_N,C_N>0$ such that, for every $u\in H^{s,p}(M)$ and every $v\in B^t_{p,p}(N)$,
\[
c_M\|u\|_{H^{s,p}(M)}
\leq
\left(\sum_{\alpha\in I}
\|(h_\alpha u)\circ\kappa_\alpha\|_{H^{s,p}(\mathbb R^n)}^p
\right)^{\frac{1}{p}}
\leq
C_M\|u\|_{H^{s,p}(M)}
\]
and
\[
c_N\|v\|_{B^t_{p,p}(N)}
\leq
\left(\sum_{i\in I_N}
\|(h_i^Nv)\circ\kappa_i^N\|_{B^t_{p,p}(\mathbb R^k)}^p
\right)^{\frac{1}{p}}
\leq
C_N\|v\|_{B^t_{p,p}(N)}.
\]
The four constants depend on $s,t,p$ and the uniform data of the pair $(M,N)$, but not on $u$, $v$, or the chart indices.

\begin{theorem}[Bessel trace on a submanifold]
\label{teo:traza-bessel-subvariedad-geometria-acotada}
\index{trace theorem!in bounded geometry}
Let \((M,N)\) be a pair of bounded geometry in the sense of Definition~\ref{def:par-geometria-acotada-traza}; in particular, \(M\) and \(N\) are smooth manifolds without boundary of dimensions \(n\) and \(k\), respectively. Let \(c=n-k\), \(1<p<\infty\), and \(s>\frac{c}{p}\). There exists a unique continuous linear operator
\[
\operatorname{Tr}_N\colon
H^{s,p}(M)
\longrightarrow
B^{s-\frac{c}{p}}_{p,p}(N)
\]
agreeing with pointwise restriction on \(C_c^\infty(M)\). It is surjective and admits a continuous linear operator
\[
\operatorname{Ex}_{M,N}\colon
B^{s-\frac{c}{p}}_{p,p}(N)
\longrightarrow
H^{s,p}(M)
\]
such that
\[
\operatorname{Tr}_N\operatorname{Ex}_{M,N}
=
\operatorname{id}_{B^{s-\frac{c}{p}}_{p,p}(N)}.
\]
\end{theorem}

\begin{proof}
Set \(t:=s-\frac{c}{p}\). For \(u\in H^{s,p}(M)\), define the localizations
\[
u_i:=\bigl(h_i u\bigr)\circ\kappa_i,
\qquad i\in I_N,
\]
regarded as distributions on \(\mathbb R^n\) by extension by zero outside the coordinate domain. This is legitimate because \(\operatorname{supp}(h_i)\Subset U_i\). By admissibility and the localization principle,
\begin{equation}
\displaystyle\sum_{i\in I_N}\|u_i\|_{H^{s,p}(\mathbb R^n)}^p
\leq
C\|u\|_{H^{s,p}(M)}^p.
\label{eq:localizacion-H-fermi}
\end{equation}
Define
\begin{equation}
\operatorname{Tr}_Nu
:=
\displaystyle\sum_{i\in I_N}
h_i^N
\left[
\bigl(\operatorname{tr}_c u_i\bigr)
\circ(\kappa_i^N)^{-1}
\right].
\label{eq:def-traza-subvariedad}
\end{equation}
The sum is locally finite. Moreover, each \(\operatorname{tr}_cu_i\) belongs to \(B^t_{p,p}(\mathbb R^k)\) by Theorem~\ref{teo:traza-euclidiana-codimension-c}.

We prove the global estimate. For \(j\in I_N\), let
\[
A(j):=\left\{i\in I_N\middle|U_i^N\cap U_j^N\neq\varnothing\right\}.
\]
The cardinality of \(A(j)\) is bounded by the multiplicity \(L\). Localizing \eqref{eq:def-traza-subvariedad} with \(h_j^N\) and transporting by \(\kappa_j^N\) to \(\mathbb R^k\) gives
\begin{equation}
\bigl(h_j^N\operatorname{Tr}_Nu\bigr)\circ\kappa_j^N
=
\displaystyle\sum_{i\in A(j)}
a_{ij}
\left[
\bigl(\operatorname{tr}_cu_i\bigr)\circ\Phi_{ij}
\right],
\label{eq:traza-local-cambio}
\end{equation}
where
\[
a_{ij}:=(h_j^Nh_i^N)\circ\kappa_j^N,
\qquad
\Phi_{ij}:=(\kappa_i^N)^{-1}\circ\kappa_j^N
\]
on the overlap, after inserting cutoffs equal to one on the supports involved. Proposition~\ref{prop:atlas-fermi-uniforme} uniformly bounds all derivatives of \(a_{ij}\), \(\Phi_{ij}\), and \(\Phi_{ij}^{-1}\). The localized version of Lemma~\ref{lema: multiplicadores y difeomorfismos sobolev}, applied on \(B^t_{p,p}(\mathbb R^{k})\), and Theorem~\ref{teo:traza-euclidiana-codimension-c} give
\begin{equation}
\left\|
\bigl(h_j^N\operatorname{Tr}_Nu\bigr)\circ\kappa_j^N
\right\|_{B^t_{p,p}(\mathbb R^k)}^p
\leq
C L^{p-1}
\displaystyle\sum_{i\in A(j)}
\|u_i\|_{H^{s,p}(\mathbb R^n)}^p.
\label{eq:cota-traza-local}
\end{equation}
Sum over \(j\). Each \(i\) appears in at most \(L\) sets \(A(j)\); by \eqref{eq:localizacion-H-fermi},
\begin{equation}
\|\operatorname{Tr}_Nu\|_{B^t_{p,p}(N)}^p
\leq
C\displaystyle\sum_{i\in I_N}\|u_i\|_{H^{s,p}(\mathbb R^n)}^p
\leq
C'\|u\|_{H^{s,p}(M)}^p.
\label{eq:cota-traza-global}
\end{equation}
Thus, \(\operatorname{Tr}_N\) is well defined and continuous. If \(u\in C_c^\infty(M)\), Lemma~\ref{lem:compatibilidades-traza-euclidiana} and \eqref{eq:particion-cuadrada-N} show that
\[
\operatorname{Tr}_Nu
=
\displaystyle\sum_{i\in I_N}(h_i^N)^2u\restriction_N
=
u\restriction_N.
\]
Density of \(C_c^\infty(M)\) in \(H^{s,p}(M)\), obtained through the trivialization and Proposition~\ref{prop:estructura-escala-bessel-euclidiana}, proves uniqueness.

We now construct a right inverse. Choose a fixed function \[
\psi\in C_c^\infty\!\left(
B_{\mathrm{euc}}^c(0,2r_F)\times B_{\mathrm{euc}}^k(0,2r_F)
\right)
\] equal to one on a common neighborhood of all supports \((h_i\circ\kappa_i)\); this is possible because these supports were obtained from the same rescaled cutoff function. For \(f\in B^t_{p,p}(N)\), set
\[
f_i:=(h_i^Nf)\circ\kappa_i^N
\]
and define
\begin{equation}
\operatorname{Ex}_{M,N}f
:=
\displaystyle\sum_{i\in I_N}
h_i
\left[
\bigl(\psi\,\operatorname{ex}_cf_i\bigr)\circ\kappa_i^{-1}
\right],
\label{eq:extension-subvariedad}
\end{equation}
extending each summand by zero outside \(U_i\). The support of \(\psi\) makes this extension smooth at the test-function level and, by duality, well defined at the distributional level.

The continuity proof repeats \eqref{eq:traza-local-cambio}--\eqref{eq:cota-traza-global} with the spaces interchanged. Uniform continuity of multipliers and coordinate changes, the multiplicity \(L\), and continuity of \(\operatorname{ex}_c\) give
\begin{equation}
\begin{split}
\|\operatorname{Ex}_{M,N}f\|_{H^{s,p}(M)}^p
&\leq
C\displaystyle\sum_{i\in I_N}
\|\psi\,\operatorname{ex}_cf_i\|_{H^{s,p}(\mathbb R^n)}^p\\
&\leq
C'\displaystyle\sum_{i\in I_N}\|f_i\|_{B^t_{p,p}(\mathbb R^k)}^p
\leq
C''\|f\|_{B^t_{p,p}(N)}^p.
\end{split}
\label{eq:cota-extension-subvariedad}
\end{equation}
Finally, \(\psi(0,u)=1\) on the relevant supports, and \(\operatorname{tr}_c\operatorname{ex}_c=\operatorname{id}\). Applying \eqref{eq:def-traza-subvariedad} to \eqref{eq:extension-subvariedad} and using \eqref{eq:particion-cuadrada-N} again gives
\[
\operatorname{Tr}_N\operatorname{Ex}_{M,N}f
=
\displaystyle\sum_{i\in I_N}(h_i^N)^2f=f.
\]
This proves both surjectivity and existence of the right inverse.
\end{proof}

\begin{theorem}[Traces in the Besov and Triebel--Lizorkin scales]
\label{teo:traza-BF-subvariedad-geometria-acotada}
Let \((M,N)\) be a pair of bounded geometry in the sense of Definition~\ref{def:par-geometria-acotada-traza}; in particular, \(M\) and \(N\) are smooth manifolds without boundary of dimensions \(n\) and \(k\), respectively. Let \(c=n-k\), \(1<p<\infty\), \(1\leq q\leq\infty\), and \(s>\frac{c}{p}\). Then there exist retractions
\[
\operatorname{Tr}_N\colon
F^s_{p,q}(M)\longrightarrow B^{s-\frac{c}{p}}_{p,p}(N),
\qquad
\operatorname{Tr}_N\colon
B^s_{p,q}(M)\longrightarrow B^{s-\frac{c}{p}}_{p,q}(N),
\]
both with continuous linear right extension operators.
\end{theorem}

\begin{proof}
Set \(t=s-\frac{c}{p}\). Retain the cutoffs and Fermi charts from the proof of Theorem~\ref{teo:traza-bessel-subvariedad-geometria-acotada}. For $u\in F^s_{p,q}(M)$, define directly, in the distributional sense, the components
\[
u_i=(h_i u)\circ\kappa_i,
\qquad i\in I_N.
\]
Formula \eqref{eq:def-traza-subvariedad} defines the trace. If \(A(j)\) is the set of indices whose charts meet chart \(j\) of \(N\), the chart-change identities of Lemma~\ref{lem:compatibilidades-traza-euclidiana}, the uniform multipliers, and the Euclidean estimate give
\[
\bigl\|(h_j^N\operatorname{Tr}_Nu)\circ\kappa_j^N
\bigr\|_{B^t_{p,p}(\mathbb R^k)}
\leq
C\displaystyle\sum_{i\in A(j)}
\|u_i\|_{F^s_{p,q}(\mathbb R^n)}.
\]
The cardinality of \(A(j)\) is at most the uniform multiplicity \(L\). Raising to the power \(p\), summing over \(j\), and observing that each \(i\) appears in at most \(L\) sets \(A(j)\) gives
\[
\|\operatorname{Tr}_Nu\|_{B^t_{p,p}(N)}^p
\leq
C\displaystyle\sum_{i\in I_N}\|u_i\|_{F^s_{p,q}(\mathbb R^n)}^p
\leq C'\|u\|_{F^s_{p,q}(M)}^p.
\]
The strong localization principle, Theorem~\ref{teo:principio-localizacion-F-H-W}, justifies the last inequality and convergence of the reconstruction.

For \(f\in B^t_{p,p}(N)\), set \(f_i=(h_i^Nf)\circ\kappa_i^N\) and use \eqref{eq:extension-subvariedad}. The synthesis estimate, uniform multiplicity, and Euclidean coretraction give
\[
\begin{split}
\|\operatorname{Ex}_{M,N}f\|_{F^s_{p,q}(M)}^p
&\leq C\displaystyle\sum_{i\in I_N}
\|\psi\operatorname{ex}_cf_i\|_{F^s_{p,q}(\mathbb R^n)}^p\\
&\leq C'\displaystyle\sum_{i\in I_N}
\|f_i\|_{B^t_{p,p}(\mathbb R^k)}^p
\leq C''\|f\|_{B^t_{p,p}(N)}^p.
\end{split}
\]
This yields the retraction and coretraction in the Triebel--Lizorkin scale.

For Besov, choose \(s_0<s<s_1\), still with \(s_0>\frac{c}{p}\), and \(s=(1-\theta)s_0+\theta s_1\). The preceding operators, constructed using the same charts, cutoffs, and Euclidean coretractions, are continuous at both endpoints. Theorem~\ref{teo:interpolacion-operadores-metodo-K} and the identities
\[
(F^{s_0}_{p,p}(M),F^{s_1}_{p,p}(M))_{\theta,q}=B^s_{p,q}(M),
\]
\[
(B^{s_0-\frac{c}{p}}_{p,p}(N),B^{s_1-\frac{c}{p}}_{p,p}(N))_{\theta,q}
=B^{s-\frac{c}{p}}_{p,q}(N)
\]
give continuity of both \(\operatorname{Tr}_N\) and \(\operatorname{Ex}_{M,N}\) in the Besov scale.

Finally, in all scales, compatibility of the trace with multipliers and \(\operatorname{tr}_c\operatorname{ex}_c=I\) allows the direct calculation
\[
\operatorname{Tr}_N\operatorname{Ex}_{M,N}f
=\displaystyle\sum_{i\in I_N}(h_i^N)^2f=f.
\]
Equality in distributions also justifies the case \(q=\infty\), where density in the space norm has not been assumed.
\end{proof}

For noninteger orders $s=m+\sigma$, with $m\in\mathbb N_0$ and $0<\sigma<1$, fix a uniform quadratic system $\mathcal Q=(\kappa_i,h_i)_{i\in I}$, with $\kappa_i$ mapping the coordinate model into $M$, and define
\[
 \|u\|_{W^{s,p}_{\mathcal Q}(M)}
 :=\left(\sum_{i\in I}
       \|(h_i u)\circ\kappa_i\|_{W^{s,p}(\mathbb R^n)}^p
   \right)^{1/p}.
\]
Localized functions are extended by zero outside their charts, and
$W^{s,p}(M)$ consists of the distributions for which this norm is finite.
Proposition~\ref{prop:slobodeckij-uniforme-localizacion} proves independence
of the system and equivalence with the intrinsic difference norm. That
comparison is proved directly from local norms and does not use the trace
corollary below.

\begin{corollary}[\(W^{s,p}(M)\) and integer-order scales]
\label{cor:traza-W-subvariedad-geometria-acotada}
Let \((M,N)\) be a pair of bounded geometry consisting of smooth manifolds without boundary, let \(c=\dim M-\dim N\), and let \(1<p<\infty\). If \(s>\frac{c}{p}\) and \(s\notin\mathbb N\), then
\[
\operatorname{Tr}_N\colon
W^{s,p}(M)\longrightarrow B^{s-\frac{c}{p}}_{p,p}(N)
\]
is a retraction. If \(m\in\mathbb N\) and \(m>\frac{c}{p}\), then
\[
\operatorname{Tr}_N\colon
W^{m,p}(M)\longrightarrow B^{m-\frac{c}{p}}_{p,p}(N)
\]
is a retraction.
\end{corollary}

\begin{proof}
For $s\notin\mathbb N$, use the global comparison $W^{s,p}(M)=B^s_{p,p}(M)$ of Proposition~\ref{prop:slobodeckij-uniforme-localizacion} and apply the Besov retraction in the preceding theorem with $q=p$. For $m\in\mathbb N$, Proposition~\ref{prop:covariantes-localizacion-orden-entero} gives $W^{m,p}(M)=H^{m,p}(M)$; then apply Theorem~\ref{teo:traza-bessel-subvariedad-geometria-acotada}.
\end{proof}

\section{Manifolds with boundary and bounded geometry}

On a manifold with boundary, geodesic balls centered at boundary points are unavailable. The appropriate local models are half-spaces, obtained through a uniform geodesic collar.

\begin{definition}[Uniform collar and uniform geodesic collar]
\label{def:collar-uniforme-geodesico}
\index{uniform collar}
Let $M$ be a smooth manifold with nonempty boundary. A \textit{uniform collar} of radius $r_\partial>0$ is a diffeomorphism
\[
\mathcal C\colon \partial M\times[0,r_\partial)\longrightarrow U_\partial,
\]
where $U_\partial$ is an open neighborhood of $\partial M$ in $M$ and $\mathcal C(x,0)=x$. If $(M,\mathbf{g})$ is Riemannian, $\boldsymbol{\nu}$ is the inward unit normal, and
\[
\mathcal C(x,t)=\exp_x(t\boldsymbol{\nu}_x),
\]
we say the collar is \textit{uniform geodesic}.
\end{definition}

\begin{definition}[Bounded geometry with boundary]
\label{def:geometria-acotada-con-frontera-rigurosa}
\index{bounded geometry!with boundary}
Let $(M,\mathbf{g})$ be a smooth Riemannian manifold of dimension $n$ with nonempty boundary, and let $\boldsymbol{\nu}$ be the inward unit normal. We say that $M$ has \textit{bounded geometry with nonempty boundary} if the following conditions hold.
\begin{enumerate}
\item There exists $r_\partial>0$ such that
\begin{equation}
\mathcal C\colon \partial M\times[0,r_\partial)\longrightarrow M,
\qquad
\mathcal C(x,t):=\exp_x(t\boldsymbol{\nu}_x),
\label{eq:collar-geodesico}
\end{equation}
is a uniform geodesic collar.
\item $\operatorname{inj}(\partial M)>0$ for the induced metric.
\item For each $0<r<r_\partial$, there exists $i(r)>0$ such that, for every
\[
x\in M\setminus\mathcal C(\partial M\times[0,r)),
\]
the mapping $\exp_x$ is defined and is a diffeomorphism from $B_{\mathbf{g}_x}(0,i(r))\subseteq T_xM$ onto its image.
\item For each $m\in\mathbb N_0$,
\[
\sup_{x\in M}|(\nabla^M)^m\mathbf{R}^M(x)|_{\mathbf{g}}<\infty.
\]
\item If $\mathbf{II}$ is the second fundamental form of $\partial M$, then, for each $m\in\mathbb N_0$,
\[
\sup_{x\in\partial M}
\bigl|(\nabla^{\mathbf{II}})^m\mathbf{II}(x)\bigr|_{\mathbf{g}}<\infty.
\]
\end{enumerate}
This is the intrinsic formulation of \cite[Definition~2.2]{Schick2001}; the derivatives in the last item are those of Definition~\ref{def:derivadas-covariantes-segunda-forma-fundamental}.
\end{definition}

\begin{proposition}[Induced geometry on the boundary]
\label{prop:frontera-geometria-acotada-inducida}
Let $(M,\mathbf{g})$ be a Riemannian manifold of bounded geometry with nonempty boundary. Then $(\partial M,\mathbf{g}\restriction_{\partial M})$ is a manifold without boundary and with bounded geometry.
\end{proposition}

\begin{proof}
The second item of the definition gives $\operatorname{inj}(\partial M)>0$. The Gauss equation \eqref{eq:gauss-traza}, applied to $N=\partial M$, expresses $\mathbf{R}^{\partial M}$ as the tangential restriction of $\mathbf{R}^M$ plus quadratic contractions of $\mathbf{II}$. Differentiating it $\ell$ times, calculation \eqref{eq:derivadas-restriccion-curvatura} expresses $(\nabla^{\partial M})^\ell \mathbf{R}^{\partial M}$ as a finite sum of contractions of $(\nabla^M)^a\mathbf{R}^M$ and $(\nabla^{\mathbf{II}})^b\mathbf{II}$, with $a,b\leq\ell$. The fourth and fifth items bound all these terms uniformly. Definition~\ref{def: geometria acotada} concludes the proof.
\end{proof}

\begin{remark}
Schick also proves that the preceding formulation is equivalent to uniform bounds on the metric and its inverse in boundary normal coordinates \cite[Theorem~2.5]{Schick2001}. Amann assumes bounded geometry already induced on $\partial M$ and proves that his formulation is equivalent to Schick's; see \cite[Definition~2.1 and Theorem~2.2]{Amann2025FunctionSpaces}. In the following proofs, we will use the five items of the definition and the preceding proposition directly.
\end{remark}

\begin{proposition}[Uniform geodesic--Fermi atlas and collar]
\label{prop:atlas-frontera-uniforme}
Let $M$ be a manifold of bounded geometry with nonempty boundary. There exist fixed radii $r_B,r_I>0$, a uniformly locally finite cover, and charts of two types:
\begin{equation}
\kappa_i^\partial\colon
B_{\mathrm{euc}}^{n-1}(0,2r_B)\times[0,2r_B)
\longrightarrow U_i^\partial,
\qquad i\in I_\partial,
\label{eq:carta-frontera-fermi}
\end{equation}
and interior geodesic charts $\kappa_j^0\colon B_{\mathrm{euc}}^n(0,2r_I)\longrightarrow U_j^0$, $j\in I_0$, such that
\[
\kappa_i^\partial(u,t)
=
\mathcal C\bigl(\exp_{p_i}^{\partial M}(\lambda_i u),t\bigr).
\]
In these charts, the uniform estimates \eqref{eq:elipticidad-fermi} hold for the metric, its inverse, and all their derivatives. There is also a family \(h_\alpha\in C_c^\infty(U_\alpha)\) such that
\[
\displaystyle\sum_{\alpha\in I}h_\alpha^2=1,
\]
and its local representations have uniformly bounded derivatives. The family \((U_\alpha,\kappa_\alpha^{-1},h_\alpha)\) is an admissible quadratic localization system, and \((U_\alpha,\kappa_\alpha^{-1},h_\alpha^2)\) is the associated admissible trivialization. Moreover, for some \(0<r_0<r_B\),
\begin{equation}
\partial_t^\ell(h_i\circ\kappa_i^\partial)(u,t)=0,
\qquad
i\in I_\partial,\quad \ell\geq1,\quad 0\leq t<r_0.
\label{eq:cortes-constantes-normal}
\end{equation}
\end{proposition}

\begin{proof}
By Proposition~\ref{prop:frontera-geometria-acotada-inducida}, $\partial M$ has bounded geometry. First consider boundary geodesic charts centered at an arbitrary point, with a common radius smaller than the uniform radius in Theorem~\ref{teo: equivalencias de geometria acotada}. The collar formula gives a Fermi chart at every center. We obtain its estimates before choosing the net of centers; thus, radius reductions will not affect any covering property already fixed.

In geodesic coordinates of \(\partial M\), set \[
x_i(u):=\exp_{p_i}^{\partial M}(\lambda_i u)
\] and consider the scaled family
\[
\widehat F_i(r,u,s)
:=
\exp^M_{x_i(u)}(sr\,\boldsymbol{\nu}_{x_i(u)}),
\qquad 0\leq s\leq1.
\]
Then \(\kappa_i^\partial(u,r)=\widehat F_i(r,u,1)\). Apply to this family the ordered-word induction \eqref{eq:palabras-ordenadas-fermi}--\eqref{eq:recurrencia-datos-fermi}, with one normal letter \(\mathbf{n}=\boldsymbol{\partial}_r\) and tangential letters \(\mathbf{t}_\mu=\boldsymbol{\partial}_{u^\mu}\). The connections are \(\nabla^M\), the Levi--Civita connection of \(\partial M\), the normal connection of \(\partial M\subset M\), and the pullback connection under \(\widehat F_i\).

To verify the sign explicitly, if
\[
\mathbf{J}_X(t):=d\mathcal C_{(x,t)}(X,0),
\qquad
\gamma_x(t)=\mathcal C(x,t),
\]
then
\[
\frac{D^2\mathbf{J}_X}{dt^2}
+\mathbf{R}^M(\mathbf{J}_X,\dot\gamma_x)\dot\gamma_x=0,
\qquad
\mathbf{J}_X(0)=X,
\qquad
\frac{D\mathbf{J}_X}{dt}(0)
=
\nabla_X^M\boldsymbol{\nu}=-\mathbf{A}_{\boldsymbol{\nu}}X.
\]
Indeed, \[
\nabla_X^M\boldsymbol{\nu}=-\mathbf{A}_{\boldsymbol{\nu}}X+\nabla_X^\perp\boldsymbol{\nu},
\qquad
\mathbf{g}(\nabla_X^M\boldsymbol{\nu},\boldsymbol{\nu})
=\frac12X\mathbf{g}(\boldsymbol{\nu},\boldsymbol{\nu})=0,
\] and the normal bundle has rank one, so \(\nabla_X^\perp\boldsymbol{\nu}=0\).

For the scaled family, if \(\mathbf{X}_\mu=dx_i(\boldsymbol{\partial}_{u^\mu})\), the length-one data are exactly
\[
\mathbf{H}_{(\mathbf{n})}=0,\qquad
\mathbf{K}_{(\mathbf{n})}=\boldsymbol{\nu},
\]
and
\[
\mathbf{H}_{(\mathbf{t}_\mu)}=\mathbf{X}_\mu,\qquad
\mathbf{K}_{(\mathbf{t}_\mu)}
=
-r\mathbf{A}_{\boldsymbol{\nu}}\mathbf{X}_\mu.
\]
The data for every longer word are
\[
\mathbf{H}_{(\mathfrak a,I)}
=
\nabla^0_{\mathfrak a}\mathbf{H}_I,
\qquad
\mathbf{K}_{(\mathfrak a,I)}
=
\nabla^0_{\mathfrak a}\mathbf{K}_I
-\mathbf{R}^M(\mathbf{H}_{(\mathfrak a)},r\boldsymbol{\nu})\mathbf{H}_I.
\]
Thus, derivatives of \(\mathbf{A}_{\boldsymbol{\nu}}\), or equivalently covariant derivatives of \(\mathbf{II}\), appear in the initial data with the sign fixed above. The interior sources are exactly
\[
\mathbf{Q}_{(\mathfrak a)}=0,\qquad
\mathbf{Q}_{(\mathfrak a,I)}
=
\nabla_{\mathfrak a}\mathbf{Q}_I+\mathscr C_{\mathfrak a}(\mathbf{Y}_I),
\]
with \(\mathscr C_{\mathfrak a}\) given by \eqref{eq:fuente-conmutacion-fermi}, which collects the commutation terms. Consequently, at order \(m\), the induction uses only \((\nabla^M)^q\mathbf{R}^M\), \(q\leq m-1\), derivatives of \(\mathbf{II}\) present in the data, and word fields of length less than \(m\).

The Grönwall-based induction of Proposition~\ref{prop:atlas-fermi-uniforme} bounds all covariant words of $\widehat F_i$. At this point, retain that covariant formulation: ambient normal coordinates centered at a point of $\partial M$ are unnecessary. Ordinary derivatives of metric coefficients are obtained by applying the Leibniz rule to inner products of variational fields, as we do next.

Throughout the collar,
\[
g_{tt}=1,\qquad g_{\mu t}=0.
\]
The first identity follows because \(t\mapsto\mathcal C(x,t)\) is a unit-speed geodesic. For the second, if \(\mathbf{J}_\mu=d\kappa_i^\partial(\boldsymbol{\partial}_{u^\mu})\), then
\[
\frac d{dt}\mathbf{g}(\mathbf{J}_\mu,\boldsymbol{\partial}_t)
=
\mathbf{g}(\nabla_{\boldsymbol{\partial}_t}\mathbf{J}_\mu,\boldsymbol{\partial}_t)
=
\mathbf{g}(\nabla_{\mathbf{J}_\mu}\boldsymbol{\partial}_t,\boldsymbol{\partial}_t)
=
\frac12\mathbf{J}_\mu\mathbf{g}(\boldsymbol{\partial}_t,\boldsymbol{\partial}_t)=0,
\]
and orthogonality holds at \(t=0\). Moreover, \(g_{\mu\nu}=\mathbf{g}(\mathbf{J}_\mu,\mathbf{J}_\nu)\), so the preceding induction controls all derivatives of the metric.

At $t=0$, the differential is the orthogonal sum of the differential of the geodesic chart of $\partial M$ and the unit normal. To measure it in different fibers, parallel transport each Jacobi field to its initial fiber. The integral equation and bounds already obtained give
\[
 \bigl|P_{t\to0}\mathbf J_X(t)-X\bigr|
 \leq Ct|X|,
\]
with $C$ independent of the center. Choose $r_B>0$ so that $8Cr_B<1/2$, $8r_B$ lies within the original collar, and all larger tangential charts lie within the uniform radius of $\partial M$. Consequently, the Fermi metric is uniformly comparable to $\mathbf g|_{\partial M}+dt^2$ in the collar of width $8r_B$. Formula $g_{\mu\nu}=\langle\mathbf J_\mu,\mathbf J_\nu\rangle$ and the Leibniz rule bound all its derivatives. Differentiating $G^{-1}G=I$ gives the corresponding inverse bounds. In particular, there exist uniform constants $c,C>0$ such that
\[
 c\,d\lambda_{\mathbf g|_{\partial M}}\,dt
 \leq d\lambda_{\mathbf g}
 \leq C\,d\lambda_{\mathbf g|_{\partial M}}\,dt.
\]
The Jacobians and their derivatives are also bounded, since the determinant is polynomial in the entries and the square root is evaluated on a compact interval in $(0,\infty)$.

The normal variable is preserved exactly under boundary chart transitions. If \(\tau_{ji}=x_j^{-1}\circ x_i\), uniqueness of collar coordinates gives
\begin{equation}
(\kappa_j^\partial)^{-1}\circ\kappa_i^\partial(u,t)
=
(\tau_{ji}(u),t).
\label{eq:transicion-collar-preserva-t}
\end{equation}
Thus, transitions between boundary charts and their inverses have uniform estimates and preserve $t$ exactly.

In choosing $r_B$, also require that all its fixed multiples used below lie within the uniform boundary radius. Then choose a maximal $r_B$-separated net in $\partial M$. Its balls of radius $r_B/2$ are disjoint, those of radius $r_B$ cover $\partial M$, and the Fermi charts cover $\mathcal C(\partial M\times[0,2r_B))$.

For the interior part, set
\[
 K:=M\setminus\mathcal C(\partial M\times[0,r_B/2)).
\]
The third item of Definition~\ref{def:geometria-acotada-con-frontera-rigurosa} provides the radius $i(r_B/4)$. Choose
\[
 0<r_I<\min\{r_B/16,\ i(r_B/4)/20\},
\]
and, if necessary, reduce it to apply the normal-coordinate estimates obtained from curvature in Chapter~\ref{cap:geometria-acotada}. Take a maximal $r_I$-separated net in $K$. The interior charts are exponentials of radius $2r_I$ centered at this net. Although $K$ has an artificial boundary, the exponentials used are those of $M$ and are defined on full tangent balls by the preceding hypothesis. The Gauss lemma and the first-exit-time argument identify their images with metric balls of that radius. Balls of radius $r_I$ cover $K$, and the larger charts do not meet $\partial M$.

We also justify separation from the smaller collar. In the controlled collar, $|dt|_{\mathbf g}=1$. A curve crossing from $t\geq r_B/2$ to $t<r_B/4$ has length at least $r_B/4$. A curve coming from outside the collar must first cross this same interval of values of $t$. It cannot leave the collar toward infinity in finite length: comparison with $\mathbf g|_{\partial M}+dt^2$ and completeness of $\partial M$ prevent this. Since $2r_I<r_B/8$, all interior charts remain outside $\mathcal C(\partial M\times[0,r_B/4))$. The union of the two families covers $M$; the band $r_B/2\leq t<r_B$ lies in regions covered by both.

Multiplicity among boundary charts is estimated on $\partial M$, since an overlap point has the same normal footpoint. Among interior charts, compare disjoint balls of radius $r_I/2$ with a ball of radius $5r_I$, within the radius $i(r_B/4)$. For mixed overlaps, an interior chart meeting a boundary chart lies in the controlled collar, and its projection onto $\partial M$ lies in a ball of radius $C r_B$. The projection is uniformly Lipschitz by the collar metric comparison. Thus, all disjoint interior balls involved lie in a half-tube of volume at most $C'r_B^n$, and each has volume at least $c'r_I^n$. Their number is bounded by $(C'/c')(r_B/r_I)^n$. In the reverse direction, the boundary centers of charts meeting an interior chart lie in a ball of $\partial M$ of radius $C''r_B$; packing on $\partial M$ bounds their number. Both nets are countable by second countability. This proves uniform local finiteness and the cross-overlap bounds.

On a mixed overlap, the connection transformation identity gives, for the transition $\Phi$ between the two charts,
\[
 \partial_a\partial_b\Phi^k
 =\Gamma_{ab}^c\partial_c\Phi^k
  -\widetilde\Gamma_{ij}^k(\Phi)
       \partial_a\Phi^i\partial_b\Phi^j.
\]
Metric comparison bounds $D\Phi$ and its inverse. Successive differentiation of this formula bounds all higher derivatives, since the metrics and their inverses have uniform jets in both systems. Order zero is controlled by their fixed-radius domains. This also gives estimates for mixed transitions, without extending a boundary chart to a normal ball centered outside $M$.

Fix $0<r_0<r_B/4$. For the partition, choose cutoffs \(\eta_i^\partial(u)\) on \(\partial M\) and a function \(\chi(t)\) equal to one for \(0\leq t<r_0\). In boundary charts use \(\eta_i^\partial(u)\chi(t)\); interior cutoffs vanish on the smaller collar. Normalizing by \[
S^{\frac{1}{2}}
=
\left(\displaystyle\sum_{\alpha\in I}\widetilde\eta_\alpha^2\right)^{\frac{1}{2}},
\] gives \(h_\alpha=\widetilde\eta_\alpha S^{-\frac{1}{2}}\). On \(0\leq t<r_0\), every summand of \(S\) contains the same factor \(\chi(t)^2=1\), and no interior cutoffs contribute; therefore, \(h_i\circ\kappa_i^\partial\) is independent of \(t\) in that region. This proves \eqref{eq:cortes-constantes-normal}. The remaining bounds follow through the quadratic normalization used in Proposition~\ref{prop:atlas-fermi-uniforme}. For \ref{item:B3}, before transporting the cutoffs choose slightly wider model functions equal to one on neighborhoods of the supports of \(\eta_i^\partial(u)\chi(t)\) and the interior cutoffs, respectively. Their transports provide functions \(\chi_\alpha\in C_c^\infty(U_\alpha)\) equal to one near \(\operatorname{supp}(h_\alpha)\); since they come from a finite family of fixed profiles, all their coordinate derivatives are bounded by constants independent of \(\alpha\).
\end{proof}

\begin{proposition}[The compact case with boundary]
\label{prop:compacta-frontera-geometria-acotada}
\index{bounded geometry!compact case with boundary}
Every compact Riemannian manifold with smooth nonempty boundary satisfies Definition~\ref{def:geometria-acotada-con-frontera-rigurosa}.
\end{proposition}

\begin{proof}
The boundary is a compact Riemannian manifold without boundary; Theorem~\ref{teo:geometria-acotada-variedad-riemanniana-compacta-geometria-acotada} gives, in particular, $\operatorname{inj}(\partial M)>0$. Apply Theorem~\ref{teo:variedades-riemannianas-existencia-de-vecindades-tubulares} to $\partial M$ in the smooth double of $M$, equipped with a smooth extension of $\mathbf{g}$. Compactness of $\partial M$ provides $r_\partial>0$ such that the normal exponential is injective on normal vectors of length less than $r_\partial$. Its restriction to the inward normal is therefore the required uniform geodesic collar.

Fix $0<r<r_\partial$ and set
$K:=M\setminus\mathcal C(\partial M\times[0,r))$. This set is compact and
contained in $\operatorname{Int}(M)$. If $K$ is empty, the third item
imposes no condition. Otherwise, for each $x\in K$, apply
Lemma~\ref{lem:vecindades-normales-convexas-locales} to the manifold without
boundary $\operatorname{Int}(M)$ and an open set
$O_x\Subset\operatorname{Int}(M)$ containing $x$. This gives an open set
$W_x\subseteq O_x$ containing $x$ and $\varepsilon_x>0$ such that, for every
$y\in W_x$, $\exp_y$ is a diffeomorphism on
$B_{\mathbf g_y}(0,2\varepsilon_x)$, with image contained in $O_x$.
Choose $x_1,\ldots,x_N\in K$ such that
$K\subseteq\bigcup_{j=1}^N W_{x_j}$ and define
\[
 i(r):=\min_{1\leq j\leq N}\varepsilon_{x_j}>0.
\]
For every $y\in K$, some $W_{x_j}$ contains $y$. The restriction of
$\exp_y$ to $B_{\mathbf g_y}(0,i(r))$ is therefore a diffeomorphism onto
its image, contained in $\operatorname{Int}(M)$. This verifies the third
item with a radius independent of the center.

For each $j\in\mathbb N_0$, the continuous functions $x\mapsto|\nabla^j\mathbf{R}^M(x)|_{\mathbf{g}}$ and $y\mapsto|\nabla^j\mathbf{II}(y)|_{\mathbf{g}}$ attain a maximum on $M$ and $\partial M$, respectively. All five items are verified.
\end{proof}

\begin{corollary}[Uniform regularity in the presence of boundary]
\label{cor:geometria-acotada-frontera-regularidad-uniforme-amann}
Let $(M,\mathbf{g})$ be a manifold of bounded geometry with nonempty boundary. The Fermi and interior atlas of Proposition~\ref{prop:atlas-frontera-uniforme}, after fixed rescaling of its domains, defines a uniformly regular Riemannian structure in the sense of Definitions~\ref{def:atlas-uniformemente-regular-amann} and \ref{def:riemanniana-uniformemente-regular-amann}. Its constants depend only on the dimension, collar width, interior and boundary injectivity radii, and bounds on $\nabla^j\mathbf{R}^M$ and $\nabla^j\mathbf{II}$ through the required order. Proposition~\ref{prop:sistema-localizacion-amann}, applied to this normalized atlas, produces an admissible quadratic localization system with the uniform bounds established there.
\end{corollary}

\begin{proof}
Use the radii $r_B,r_I$ of the proposition. Choose $b_B,b_I>0$ such that $b_B\sqrt n<r_B$ and $b_I\sqrt n<r_I$. On the boundary and in the interior, respectively, define
\[
 \bar\kappa_i^{-1}(x)
 :=\kappa_i^\partial\bigl(b_B(x_2,\ldots,x_n),2r_Bx_1\bigr),
 \qquad x\in Q_+^n,
\]
\[
 \bar\kappa_j^{-1}(x):=\kappa_j^0(b_Ix),
 \qquad x\in Q^n.
\]
These are two fixed linear transformations with bounded inverses; thus, they preserve all uniform bounds, although the constants may depend on $b_B,b_I,r_B$.

Restriction to cubes is accompanied by a new choice of centers. On $\partial M$, choose a net of spacing $\delta_B$ small enough that the images of tangential cubes of half the side length cover $\partial M$. This follows from the first-exit-time argument: the metric lower bound gives a uniform minimum length for leaving the smaller cube. The boundary half-boxes thus cover $0\leq t<r_B$. On $K=M\setminus\mathcal C(\partial M\times[0,r_B/2))$, choose a net of spacing $\delta_I$ less than the minimum length required to leave an interior cube of half the side length. Its half-sized boxes cover $K$. The interior radii remain less than $r_I$, so the charts do not meet $t=0$. These two families cover $M$ after shrinking by the same factor $1/2$. The volume arguments of the proposition, applied at the new fixed scales, preserve uniform multiplicity.

The proof of the cited proposition gives, for each $j\in\mathbb N_0$, a constant $C_{\mathrm{tr}}(j)>0$, independent of the charts, bounding all coordinate transitions in $BUC^j$. It also gives constants $c_g,C_g>0$ such that
\[
 c_g|\xi|^2\leq g_\kappa(x)(\xi,\xi)
 \leq C_g|\xi|^2
\]
for every chart, every point of its normalized domain, and every $\xi\in\mathbb R^n$. For each $j$, it also provides constants $G_j,G'_j>0$, independent of the chart, satisfying
\[
 Q_\alpha:=\bar\kappa_\alpha(U_\alpha)
 =\begin{cases}
 Q_+^n,&\alpha\in I_\partial,\\
 Q^n,&\alpha\in I_0,
 \end{cases}
 \qquad I:=I_\partial\sqcup I_0,
\]
and
\[
 \|(g_{\bar\kappa_\alpha})_{ab}\|_{BUC^j(Q_\alpha)}\leq G_j,
 \qquad
 \|(g_{\bar\kappa_\alpha}^{-1})^{ab}\|_{BUC^j(Q_\alpha)}\leq G'_j,
 \qquad \alpha\in I.
\]
Uniform continuity of derivatives of order $j$ follows by applying the mean value theorem to the bounds of order $j+1$ in the slightly larger charts used in that proof. Now apply Proposition~\ref{prop:sistema-localizacion-amann} to the normalized atlas. This gives cutoffs $\pi_\alpha\in C_c^\infty(U_\alpha,[0,1])$ and a fixed coordinate cutoff $\widehat\chi$ such that
\[
 \sum_{\alpha\in I}\pi_\alpha^2=1,
 \qquad
 \widehat\chi=1\quad\hbox{on }
 \operatorname{supp}(\pi_\alpha\circ\bar\kappa_\alpha^{-1}).
\]
For each $j\in\mathbb N_0$, there exists $D_j>0$, independent of $\alpha$, with
\[
 \|\pi_\alpha\circ\bar\kappa_\alpha^{-1}\|_{BUC^j(Q_\alpha)}
 \leq D_j,
 \qquad \alpha\in I.
\]
The constant $D_j$ depends only on $j$, the multiplicity $N$, the shrinking factor $\rho=\frac{1}{2}$, and $C_{\mathrm{tr}}(0),\ldots,C_{\mathrm{tr}}(j)$. These cutoffs are constructed after reselecting the centers and restricting to cubes, so they may differ from the cutoffs $h_\beta$ in the quadratic localization system of Proposition~\ref{prop:atlas-frontera-uniforme}.
\end{proof}

\begin{theorem}[Full comparison with boundary]
\label{teo:amann-ur-equivale-geometria-acotada-frontera}
A Riemannian manifold with smooth nonempty boundary is uniformly regular if and only if it has bounded geometry with boundary in the sense of Definition~\ref{def:geometria-acotada-con-frontera-rigurosa}.
\end{theorem}

The implication from bounded geometry to uniform regularity is Corollary~\ref{cor:geometria-acotada-frontera-regularidad-uniforme-amann}, whose proof uses the Fermi atlas already constructed. Conversely, the results in \cite[Theorems~2.2 and~3.5, \S\S2--3]{Amann2025FunctionSpaces} state that a uniformly regular Riemannian atlas recovers a uniform geodesic collar, interior and boundary injectivity radii, and curvature bounds, and that this formulation is equivalent to Schick's through the second-fundamental-form bounds. The intrinsic equivalence with the second-fundamental-form formulation is also found in \cite[\S2, Definition~2.2 and Theorem~2.5]{Schick2001}.

\begin{definition}[Localized spaces on a manifold with boundary]
\label{def:espacios-geometria-acotada-con-frontera}
Let $M$ be a Riemannian manifold of bounded geometry with nonempty boundary in the sense of Definition~\ref{def:geometria-acotada-con-frontera-rigurosa}, let $\mathbb R^n_+:=\mathbb R^{n-1}\times[0,\infty)$, and fix the quadratic localization system of Proposition~\ref{prop:atlas-frontera-uniforme}.
\begin{enumerate}[label=(\alph*)]
\item For $s\in\mathbb R$ and $1<p<\infty$, define
\[
H^{s,p}(\mathbb R^n_+)
:=
\left\{U\restriction_{\mathbb R^n_+}\middle|U\in H^{s,p}(\mathbb R^n)\right\},
\]
with quotient norm
\begin{equation}
\label{eq:norma-cociente-H-semiespacio}
\|u\|_{H^{s,p}(\mathbb R^n_+)}
:=
\inf_{\substack{U\in H^{s,p}(\mathbb R^n)\\U\restriction_{\mathbb R^n_+}=u}}\|U\|_{H^{s,p}(\mathbb R^n)}.
\end{equation}
Here and below, distributional restriction is taken to the open set $\mathbb R^{n-1}\times(0,\infty)$; we retain the notation $\mathbb R^n_+=\mathbb R^{n-1}\times[0,\infty)$ because trace values on $t=0$ will be considered later. The restriction mapping is a continuous surjection, and this norm is, by definition, the quotient norm modulo its kernel. In particular, it does not depend on a chosen extension, and $H^{s,p}(\mathbb R^n_+)$ is a Banach space.

In all restriction scales of this chapter, a distribution on $M$ means a distribution on $\operatorname{Int}(M)$. No independent distributions supported on $\partial M$ are added: boundary values will be obtained through trace operators. The space $H^{s,p}(M)$ consists of those distributions $u$ for which the norm
\begin{equation}
\|u\|_{H^{s,p}(M)}
:=
\left(
\sum_{i\in I_\partial}
\|(h_i u)\circ\kappa_i^\partial\|_{H^{s,p}(\mathbb R^n_+)}^p
+
\sum_{j\in I_0}
\|(h_j u)\circ\kappa_j^0\|_{H^{s,p}(\mathbb R^n)}^p
\right)^{\frac{1}{p}}.
\label{eq:norma-H-frontera}
\end{equation}
is finite.

\item If $s\in\mathbb R$, $1\leq p<\infty$, and $1\leq q\leq\infty$, define
\[
F^s_{p,q}(\mathbb R^n_+)
:=
\left\{U\restriction_{\mathbb R^n_+}\middle|
U\in F^s_{p,q}(\mathbb R^n)\right\}
\]
and
\begin{equation}
\label{eq:norma-cociente-F-semiespacio}
\|u\|_{F^s_{p,q}(\mathbb R^n_+)}
:=
\inf_{\substack{U\in F^s_{p,q}(\mathbb R^n)\\U\restriction_{\mathbb R^n_+}=u}}\|U\|_{F^s_{p,q}(\mathbb R^n)}.
\end{equation}
Again, this is the quotient norm of the restriction operator and does not depend on a particular extension. The space $F^s_{p,q}(M)$ consists of the distributions for which
\[
\left(
\sum_{i\in I_\partial}
\|(h_i u)\circ\kappa_i^\partial\|_{F^s_{p,q}(\mathbb R^n_+)}^p
+
\sum_{j\in I_0}
\|(h_j u)\circ\kappa_j^0\|_{F^s_{p,q}(\mathbb R^n)}^p
\right)^{\frac{1}{p}}
\]
is finite; this expression is its norm.

\item For $s\in\mathbb R$, $1\leq p<\infty$, and $1\leq q\leq\infty$, also define the Euclidean restriction space
\[
B^s_{p,q}(\mathbb R^n_+)
:=
\left\{U\restriction_{\mathbb R^n_+}\middle|
U\in B^s_{p,q}(\mathbb R^n)\right\}
\]
with quotient norm
\begin{equation}
\label{eq:norma-cociente-B-semiespacio}
\|u\|_{B^s_{p,q}(\mathbb R^n_+)}
:=
\inf_{\substack{U\in B^s_{p,q}(\mathbb R^n)\\U\restriction_{\mathbb R^n_+}=u}}\|U\|_{B^s_{p,q}(\mathbb R^n)}.
\end{equation}
Provisionally fix $s_0<s<s_1$ and $\theta=\frac{s-s_0}{s_1-s_0}$. Define
\begin{equation}
\label{eq:besov-frontera-extremos-provisionales}
B^s_{p,q}(M;s_0,s_1)
:=
\bigl(F^{s_0}_{p,p}(M),F^{s_1}_{p,p}(M)\bigr)_{\theta,q}.
\end{equation}
Independence of $s_0,s_1$, which is not assumed in this definition, is proved in Corollary~\ref{cor:independencia-extremos-besov-frontera} using a common analysis--synthesis pair. After that corollary, the endpoints will be omitted and the space written simply as $B^s_{p,q}(M)$.
\end{enumerate}
Compositions with parametrizations are understood in the sense of Proposition~\ref{prop:representacion-local-distribucion-coordenadas}; in boundary charts, use the restriction norm just defined. On $\partial M$, use the definitions without boundary from the preceding chapter.
\end{definition}

\begin{lemma}[Localized operators between boundary models]
\label{lem:operadores-localizados-semiespacio-fermi}
Let $s\in\mathbb R$ and choose $k\in\mathbb N$ with $k>|s|+n+1$. Consider $H^{s,p}(\mathbb R^n_+)$ for $1<p<\infty$, $F^s_{p,q}(\mathbb R^n_+)$ for $1\leq p<\infty$, $1\leq q\leq\infty$, and $B^s_{p,q}(\mathbb R^n_+)$ for $1\leq p<\infty$, $1\leq q\leq\infty$.
\begin{enumerate}[label=(\alph*)]
\item Rychkov's cone formula defines a single linear operator
\[
 \mathbf{E}_+\colon \mathcal A(\mathbb R^n_+)\longrightarrow
 \mathcal A(\mathbb R^n),
 \qquad \mathcal A\in\{H^{s,p},F^s_{p,q},B^s_{p,q}\},
\]
such that $(\mathbf{E}_+u)\restriction_{\mathbb R^n_+}=u$. There exist constants $C_+^H(s,p)$, $C_+^F(s,p,q)$, and $C_+^B(s,p,q)$, depending additionally only on $n$ and the fixed kernels, bounding the three norms of $\mathbf{E}_+$.

\item Let $U',V'\subseteq\mathbb R^{n-1}$ be open, $\vartheta\colon U'\longrightarrow V'$ a diffeomorphism of class $C^{k+1}$, $\lambda>0$, and $\Psi(x',t):=(\vartheta(x'),\lambda t)$. Suppose there exist $\lambda_-,\lambda_+>0$ such that
\[
 \begin{gathered}
 \lambda_-\leq\lambda\leq\lambda_+,
 \qquad c_J\leq|\det D\vartheta(x')|\leq C_J
 \quad(x'\in U'),\\
 \|D^\ell\vartheta\|_{L^\infty(U')}+
 \|D^\ell\vartheta^{-1}\|_{L^\infty(V')}\leq A,
 \qquad \ell\in\{1,\ldots,k+1\}.
 \end{gathered}
\]
and $\eta\in C_c^{k+1}(U'\times[0,\infty))$ satisfies $\|\eta\|_{BC^{k+1}(U'\times[0,\infty))}\leq A_\eta$. Let $K_\eta$ be the tangential projection of $\operatorname{supp}\eta$. Also fix $R,\delta>0$ such that
\[
 \begin{gathered}
 \operatorname{supp}\eta
 \subseteq\overline B_{\mathrm{euc}}^{n-1}(0,R)\times[0,R],\\
 \operatorname{dist}(K_\eta,\mathbb R^{n-1}\setminus U')\geq\delta,\\
 \operatorname{dist}(\vartheta(K_\eta),\mathbb R^{n-1}\setminus V')\geq\delta.
 \end{gathered}
\]
These distances concern the artificial chart boundaries; the support may meet $t=0$. The initial formula
\[
 T_{\eta,\Psi}^{+\leftarrow+}u
 :=\begin{cases}
 \eta\,(u\circ\Psi),&\text{on }U'\times[0,\infty),\\
 0,&\text{elsewhere in }\mathbb R^n_+,
 \end{cases}
\]
defines continuous operators
\[
 T_{\eta,\Psi}^{+\leftarrow+}\colon
 \mathcal A(\mathbb R^n_+)\longrightarrow \mathcal A(\mathbb R^n_+)
\]
in all three scales. Their norms are bounded by constants $M_H,M_F(q),M_B(q)$ depending on $n,s,p,k,c_J,C_J,A,A_\eta,\lambda_-,\lambda_+,R,\delta$ and, in the last two scales, also on $q$.

\item Let $\mathbb X,\mathbb Y\in\{\mathbb R^n,\mathbb R^n_+\}$, let $U\subseteq\mathbb Y$ and $V\subseteq\mathbb X$ be relatively open, let $\Phi\colon U\longrightarrow V$ be a diffeomorphism, and let $\zeta\in C_c^{k+1}(U)$. If a half-space occurs and the relevant support meets $t=0$, assume $\mathbb X=\mathbb Y=\mathbb R^n_+$ and that $(\zeta,\Phi)$ has the form in (b). In the other cases, assume the relevant supports have distance at least $d>0$ from every boundary present: $\operatorname{dist}(\operatorname{supp}\zeta,\partial\mathbb R^n_+)\geq d$ if $\mathbb Y=\mathbb R^n_+$ and $\operatorname{dist}(\Phi(\operatorname{supp}\zeta),
\partial\mathbb R^n_+)\geq d$ if $\mathbb X=\mathbb R^n_+$. Then
\[
 T_{\zeta,\Phi}^{\mathbb Y\leftarrow\mathbb X}\colon
 \mathcal A(\mathbb X)\longrightarrow \mathcal A(\mathbb Y)
\]
is continuous. For uniform bounds, include the diameter of the supports and their positive separation from the relative boundaries of $U$ and $V$, together with the coordinate-change bounds required in the localized version of Lemma~\ref{lema: multiplicadores y difeomorfismos sobolev}. If $\widetilde T$ denotes the localized full-space operator associated with the same data, $r_+$ the restriction, and $\mathbf{E}_+$ the operator in (a), the four types are
\[
\begin{array}{c|c}
(\mathbb X,\mathbb Y)&T_{\zeta,\Phi}^{\mathbb Y\leftarrow\mathbb X}\\ \hline
(\mathbb R^n,\mathbb R^n)&\widetilde T\\
(\mathbb R^n_+,\mathbb R^n)&\widetilde T \mathbf{E}_+\\
(\mathbb R^n,\mathbb R^n_+)&r_+\widetilde T\\
(\mathbb R^n_+,\mathbb R^n_+)&r_+\widetilde T \mathbf{E}_+.
\end{array}
\]
In each case, on smooth functions the formula means $T_{\zeta,\Phi}^{\mathbb Y\leftarrow\mathbb X}u
=\zeta(u\circ\Phi)$ on $U$ and zero outside $U$. The preceding constants also depend on $d$ and bounds for $\Phi$, $\Phi^{-1}$, their Jacobians, and their derivatives through order $k+1$, but not on $u$.
\end{enumerate}
All these assertions include $s<0$. At endpoint fine indices, the distributional formula specifies the operator.
\end{lemma}

\begin{proof}
Specialize the construction of Lemma~\ref{lem:nucleos-rychkov-cono} to $\omega=0$. With its cone kernels, set
\[
 \mathbf{E}_+u:=\sum_{j\in\mathbb N_0}
 \psi_j*\bigl(\mathbf1_{\mathbb R^n_+}(\varphi_j*u)\bigr).
\]
The identity \eqref{eq:reproduccion-conica-rychkov} proves $r_+\mathbf{E}_+=I$. Estimates \eqref{eq:casi-ortogonalidad-conica-rychkov} and \eqref{eq:analisis-medias-locales-conicas} give $C_+^F(s,p,q)$; for $q=2$ and $H^{s,p}(\mathbb R^n)=F^s_{p,2}(\mathbb R^n)$, they give $C_+^H(s,p)$. For Besov, take the $\ell^q(\mathbb N_0)$ norm outside $L^p(\mathbb R^n)$ in the same estimates, giving $C_+^B(s,p,q)$. Fix the Rychkov kernels once, with all vanishing moments required at positive levels. To estimate a particular order $s$, only finitely many $N>|s|+1$ of these moments are used; the kernels do not change as $s$ varies. This gives a single operator at every order, including $s<0$.

We record the density that will be used. Let $Q_{\mathcal A}$ be the restriction quotient mapping from $\mathcal A(\mathbb R^n)$ to $\mathcal A(\mathbb R^n_+)$. If $\mathcal S(\mathbb R^n)$ is dense in $\mathcal A(\mathbb R^n)$ and $u=Q_{\mathcal A}U$, choose $U_j\in\mathcal S$ with $U_j\longrightarrow U$; since $\|Q_{\mathcal A}\|\leq1$, we have $Q_{\mathcal A}U_j\longrightarrow u$. Thus, restrictions of Schwartz functions are dense in $H^{s,p}(\mathbb R^n)$, in $F^s_{p,q}(\mathbb R^n)$ when $q<\infty$, and in $B^s_{p,q}(\mathbb R^n)$ when $p,q<\infty$. This conclusion is neither obtained nor needed at the other indices.

For (b), let $a_1,\ldots,a_{k+2}$ be the unique coefficients satisfying
\[
 \sum_{\nu=1}^{k+2}a_\nu(-\nu)^j=1,
 \qquad 0\leq j\leq k+1,
\]
and define
\[
 \widetilde\eta(x',t):=
 \begin{cases}
 \eta(x',t),&t\geq0,\\
 \displaystyle\sum_{\nu=1}^{k+2}a_\nu\eta(x',-\nu t),&t<0.
 \end{cases}
\]
The Vandermonde identities match derivatives through order $k+1$. Moreover,
$\widetilde\eta\in C_c^{k+1}(U'\times\mathbb R)$
and, with
\[
 C_k:=1+\sum_{\nu=1}^{k+2}|a_\nu|\max\{1,\nu^{k+1}\},
\]
we have $\|\widetilde\eta\|_{BC^{k+1}(U'\times\mathbb R)}\leq C_kA_\eta$. The diffeomorphism $\widetilde\Psi(x',t):=(\vartheta(x'),\lambda t)$ preserves the half-spaces. Moreover,
\[
 D\widetilde\Psi=
 \begin{pmatrix}D\vartheta&0\\0&\lambda\end{pmatrix},
 \qquad
 D\widetilde\Psi^{-1}=
 \begin{pmatrix}D\vartheta^{-1}&0\\0&\lambda^{-1}\end{pmatrix},
 \qquad
 |\det D\widetilde\Psi|=\lambda|\det D\vartheta|.
\]
Consequently, derivatives of $\widetilde\Psi$ and $\widetilde\Psi^{-1}$ through order $k+1$, and upper and lower bounds for their Jacobians, are controlled by constants depending only on $c_J,C_J,A,\lambda_-,\lambda_+$. The support of $\widetilde\eta$ remains in a compact set of diameter controlled by $R$ and separated from the boundary of $U'\times\mathbb R$ by $\delta$. The corresponding image margin is also uniform. The localized version of Lemma~\ref{lema: multiplicadores y difeomorfismos sobolev}, with these distances and bounds, and item (a) give
\begin{align*}
 \|r_+T_{\widetilde\eta,\widetilde\Psi}\mathbf{E}_+u\|_{H^{s,p}(\mathbb R^n_+)}
 &\leq C^H_{\mathrm{loc}}C_+^H(s,p)
 \|u\|_{H^{s,p}(\mathbb R^n_+)},\\
 \|r_+T_{\widetilde\eta,\widetilde\Psi}\mathbf{E}_+u\|_{F^s_{p,q}(\mathbb R^n_+)}
 &\leq C^F_{\mathrm{loc}}(q)C_+^F(s,p,q)
 \|u\|_{F^s_{p,q}(\mathbb R^n_+)},\\
 \|r_+T_{\widetilde\eta,\widetilde\Psi}\mathbf{E}_+u\|_{B^s_{p,q}(\mathbb R^n_+)}
 &\leq C^B_{\mathrm{loc}}(q)C_+^B(s,p,q)
 \|u\|_{B^s_{p,q}(\mathbb R^n_+)}.
\end{align*}
The restriction $r_+$ has norm at most one by the quotient definition. Thus, we may take
\[
 M_H:=C^H_{\mathrm{loc}}C_+^H(s,p),
 \quad
 M_F(q):=C^F_{\mathrm{loc}}(q)C_+^F(s,p,q),
 \quad
 M_B(q):=C^B_{\mathrm{loc}}(q)C_+^B(s,p,q).
\]
The third bound follows directly from the almost-diagonal Besov estimate in Lemma~\ref{lema: multiplicadores y difeomorfismos sobolev}; the outer sum $\ell^p(I)$ and integration in the $K$ method are not interchanged. These formulas display every stated dependence. In ranges where restrictions of Schwartz functions are dense, they also prove uniqueness of the extension; at the endpoints, the composition $r_+\widetilde T \mathbf{E}_+$ is its distributional definition.

For (c), if the supports are separated from $t=0$, the effective open sets are open in $\mathbb R^n$ and the localized operator $\widetilde T$ from the full-space lemma applies without crossing a boundary. Precomposing with $\mathbf{E}_+$ when $\mathbb X=\mathbb R^n_+$ and postcomposing with $r_+$ when $\mathbb Y=\mathbb R^n_+$ gives exactly the four factorizations in the table. If the support meets $t=0$, the construction in (b) provides $\widetilde T$. In $H$ and $F$, the norms are the product of the localized constant with $C_+^H$ or $C_+^F$ whenever $\mathbf{E}_+$ occurs. In $B$, the same procedure gives the localized constant $C^B_{\mathrm{loc}}(q)$, multiplied by $C_+^B(s,p,q)$ precisely in the two types whose domain is $\mathbb R^n_+$; the restriction $r_+$ has norm at most one. The almost-diagonal argument of the full-space lemma includes $s<0$, proving the assertions for negative orders through the same local estimates.
\end{proof}

\begin{proposition}[Independence of the boundary atlas]
\label{prop:independencia-atlas-frontera}
Let $M$ be a Riemannian manifold of bounded geometry with nonempty boundary, and let $\mathcal T_1$ and $\mathcal T_2$ be two atlases satisfying Proposition~\ref{prop:atlas-frontera-uniforme}. For each space in Definition~\ref{def:espacios-geometria-acotada-con-frontera}, there exist constants $c_{12}^{\mathcal A},C_{12}^{\mathcal A}>0$, one pair for each function-space scale $\mathcal A$, depending on its indices and the uniform data of both atlases but not on $u$, such that
\[
c_{12}^{\mathcal A}\|u\|_{\mathcal A,\mathcal T_1}
\leq \|u\|_{\mathcal A,\mathcal T_2}
\leq C_{12}^{\mathcal A}\|u\|_{\mathcal A,\mathcal T_1}.
\]
For $B^s_{p,q}(M)$, provisionally fix the same interpolation endpoints in both atlases; their independence is proved in Corollary~\ref{cor:independencia-extremos-besov-frontera}.
\end{proposition}

The retraction--coretraction argument corresponds to \cite[Theorems~4.2 and 4.4]{Amann2025FunctionSpaces}.

\begin{proof}
Let \(\mathcal T_1\) and \(\mathcal T_2\) be the two systems. Insert the identity \(\displaystyle \sum_{\beta\in I_2}h_{2,\beta}^2=1\) into each localization \(h_{1,\alpha}u\). The two cross-overlap bounds leave at most $L$ summands per index, where $L$ depends only on the uniform data of both atlases.

Between two boundary charts, using the same collar normal parameter makes the Fermi transition
\[
 (x',t)\longmapsto\bigl(\vartheta_{\alpha\beta}(x'),t\bigr).
\]
The bounds in Proposition~\ref{prop:atlas-frontera-uniforme} provide numbers $c_J,C_J,A,A_\eta>0$, independent of $\alpha,\beta$, verifying the hypotheses of Lemma~\ref{lem:operadores-localizados-semiespacio-fermi}(b). Thus, each boundary--boundary summand has norm at most $M_H$, $M_F(q)$, or $M_B(q)$ times the source localization norm, according to the scale; these assertions include $s<0$. For two interior charts, use Lemma~\ref{lema: multiplicadores y difeomorfismos sobolev} on $\mathbb R^n$. On a mixed overlap, the interior cutoff supports are separated from $t=0$ by a uniform distance; item (c) of the preceding lemma again reduces the estimate to the full-space case. Denote by $M_{\mathcal A}$ the largest constant for scale $\mathcal A$. It depends on $n,s,p$, also on $q$ in $F$ and $B$, and on the uniform data of $\mathcal T_1,\mathcal T_2$, but not on the indices.

The inequality $\displaystyle \left(\displaystyle\sum_{\ell=1}^L a_\ell\right)^p
\leq L^{p-1}\displaystyle\sum_{\ell=1}^L a_\ell^p$ and a second overlap count give
\[
\|u\|_{H^{s,p}_{\mathcal T_1}(M)}
\leq
LM_H\|u\|_{H^{s,p}_{\mathcal T_2}(M)}.
\]
Interchanging the atlases gives $\|u\|_{H^{s,p}_{\mathcal T_2}(M)}
\leq L'M'_H\|u\|_{H^{s,p}_{\mathcal T_1}(M)}$; taking $c^H_{12}:=(LM_H)^{-1}$ and $C^H_{12}:=L'M'_H$ yields the assertion in $H$. In $F$, the same count gives explicitly
\[
 c^F_{12}(q):=(LM_F(q))^{-1},
 \qquad C^F_{12}(q):=L'M'_F(q).
\]
For $B$, fix provisional endpoints $a<s<b$ and $\theta=\frac{s-a}{b-a}$. If
\[
 c_t\|u\|_{F^t_{p,p,\mathcal T_1}(M)}
 \leq\|u\|_{F^t_{p,p,\mathcal T_2}(M)}
 \leq C_t\|u\|_{F^t_{p,p,\mathcal T_1}(M)},
 \qquad t\in\{a,b\},
\]
are the inequalities just proved with $c_t:=(L_tM_F(t,p))^{-1}$ and $C_t:=L'_tM'_F(t,p)$, Theorem~\ref{teo:interpolacion-operadores-metodo-K} gives
\[
 c_a^{1-\theta}c_b^\theta
 \|u\|_{B^s_{p,q}(M;a,b),\mathcal T_1}
 \leq\|u\|_{B^s_{p,q}(M;a,b),\mathcal T_2}
 \leq C_a^{1-\theta}C_b^\theta
 \|u\|_{B^s_{p,q}(M;a,b),\mathcal T_1}.
\]
Thus, for these endpoints one may take $c^B_{12}:=c_a^{1-\theta}c_b^\theta$ and $C^B_{12}:=C_a^{1-\theta}C_b^\theta$. Subsequent independence of $a,b$ allows them to be omitted.
\end{proof}

\begin{corollary}[Functional connection with Amann's theory in the presence of boundary]
\label{cor:puente-funcional-amann-sobolev-frontera}
Let $M$ be a Riemannian manifold of bounded geometry with nonempty boundary, let $\mathfrak K=(U_\alpha,\phi_\alpha)_{\alpha\in I}$ be the normalized atlas constructed in Corollary~\ref{cor:geometria-acotada-frontera-regularidad-uniforme-amann}, and let $(\pi_\alpha)_\alpha$ be the new quadratic system obtained by applying Proposition~\ref{prop:sistema-localizacion-amann} to $\mathfrak K$. Set \[
 X_\alpha:=
 \begin{cases}
  \mathbb R^n,&U_\alpha\cap\partial M=\varnothing,\\
 \mathbb R^n_+,&U_\alpha\cap\partial M\neq\varnothing,
 \end{cases}
\] and, at a boundary index, postcompose the chart with the linear isometry \[
 \tau(x_1,\ldots,x_n):=(x_2,\ldots,x_n,x_1),
\] identifying the model $Q_+^n$ with normal variable first with the half-space $\mathbb R^n_+=\mathbb R^{n-1}\times[0,\infty)$ used in Definition~\ref{def:espacios-geometria-acotada-con-frontera}. Retain the symbol $\phi_\alpha$ for the postcomposed chart and write $\Omega_\alpha:=\phi_\alpha(U_\alpha)\subseteq X_\alpha$. For $s\in\mathbb R$ and $1<p<\infty$, define
\[
 \mathbb H^{s,p}_{\partial}(\mathfrak K)
 :=\left\{v=(v_\alpha)_{\alpha\in I}\middle|
 \sum_{\alpha\in I}\|v_\alpha\|_{H^{s,p}(X_\alpha)}^p<+\infty\right\}
\]
with norm given by the $p$th root of the preceding sum. For $s\in\mathbb R$, $1\leq p<\infty$, and $1\leq q\leq\infty$, likewise set
\[
 \mathbb F^s_{p,q,\partial}(\mathfrak K)
 :=\left\{v=(v_\alpha)_{\alpha\in I}\middle|
 \sum_{\alpha\in I}\|v_\alpha\|_{F^s_{p,q}(X_\alpha)}^p<+\infty\right\}.
\]
Both spaces carry the norm given by the $p$th root of the corresponding sum. Let $\chi\in C_c^\infty(Q^n,[0,1])$ be the fixed cutoff chosen in Proposition~\ref{prop:sistema-localizacion-amann}. Define
\[
 \widehat\chi_\alpha:=
 \begin{cases}
  \chi\restriction_{\Omega_\alpha},
   &X_\alpha=\mathbb R^n,\\
  (\chi\circ\tau^{-1})\restriction_{\Omega_\alpha},
   &X_\alpha=\mathbb R^n_+.
 \end{cases}
\]
Then $\widehat\chi_\alpha$ equals one on a neighborhood of $\operatorname{supp}(\pi_\alpha\circ\phi_\alpha^{-1})$ and, for each $j\in\mathbb N_0$,
\[
 \sup_{\alpha\in I}\|\widehat\chi_\alpha\|_{BUC^j(\Omega_\alpha)}
 \leq \|\chi\|_{BUC^j(Q^n)}<+\infty,
\]
because $\tau$ is a linear isometry. Let $\widetilde\chi_\alpha$ be the extension by zero of $\widehat\chi_\alpha$ from $\Omega_\alpha$ to $X_\alpha$. Since its support stays separated from the nongeometric faces of $\Omega_\alpha$, this extension belongs to $C_c^\infty(X_\alpha)$, including when $X_\alpha=\mathbb R^n_+$ and the support meets $t=0$. For a component $v_\alpha$ of either sequence space, define $\mathcal E_\alpha(\widetilde\chi_\alpha v_\alpha)$ as the distribution whose coordinate representation on $U_\alpha$ is $(\widetilde\chi_\alpha v_\alpha)\restriction_{\Omega_\alpha}$ and which is extended by zero outside $U_\alpha$. Likewise, the tilde in $\widetilde{(\pi_\alpha u)\circ\phi_\alpha^{-1}}$ denotes extension by zero from $\Omega_\alpha$ to $X_\alpha$. This is legitimate because $\operatorname{supp}(\pi_\alpha\circ\phi_\alpha^{-1})\Subset\Omega_\alpha$. The operators, whose formulas do not depend on $s,p,q$,
\[
 \mathcal S_{\partial,\mathfrak K}u
 :=\left(
 \widetilde{(\pi_\alpha u)\circ\phi_\alpha^{-1}}
 \right)_{\alpha\in I},
 \qquad
 \mathcal R_{\partial,\mathfrak K}v
 :=\sum_{\alpha\in I}
 \pi_\alpha\,\mathcal E_\alpha(\widetilde\chi_\alpha v_\alpha)
\]
are continuous in both the $H$ and $F$ scales, and \[
 \mathcal R_{\partial,\mathfrak K}
 \mathcal S_{\partial,\mathfrak K}=\operatorname{id}
\] on each of these spaces. There exist constants $c^\partial_{s,p},C^\partial_{s,p},D^\partial_{s,p}>0$ such that
\begin{align}
 c^\partial_{s,p}\|u\|_{H^{s,p}(M)}
 &\leq
 \|\mathcal S_{\partial,\mathfrak K}u\|_{
   \mathbb H^{s,p}_{\partial}(\mathfrak K)}
 \leq C^\partial_{s,p}\|u\|_{H^{s,p}(M)},
 \label{eq:comparacion-amann-frontera-analisis}\\
 \|\mathcal R_{\partial,\mathfrak K}v\|_{H^{s,p}(M)}
 &\leq D^\partial_{s,p}
 \|v\|_{\mathbb H^{s,p}_{\partial}(\mathfrak K)}.
 \label{eq:comparacion-amann-frontera-sintesis}
\end{align}
Moreover, there exist $c^F_{s,p,q},C^F_{s,p,q},D^F_{s,p,q}>0$ such that
\begin{align}
 c^F_{s,p,q}\|u\|_{F^s_{p,q}(M)}
 &\leq
 \|\mathcal S_{\partial,\mathfrak K}u\|_{
 \mathbb F^s_{p,q,\partial}(\mathfrak K)}
 \leq C^F_{s,p,q}\|u\|_{F^s_{p,q}(M)},
 \label{eq:comparacion-amann-frontera-analisis-F}\\
 \|\mathcal R_{\partial,\mathfrak K}v\|_{F^s_{p,q}(M)}
 &\leq D^F_{s,p,q}
 \|v\|_{\mathbb F^s_{p,q,\partial}(\mathfrak K)}.
 \label{eq:comparacion-amann-frontera-sintesis-F}
\end{align}
For a fixed integer $k>|s|+n+1$, the constants depend only on $n,s,p,k$, the two overlap bounds, $q$ in scale $F$, and uniform coordinate-transition and cutoff bounds through order $k+1$; they do not depend on $\alpha$, $u$, or $v$.
\end{corollary}

\begin{proof}
Denote by $\mathcal T_F=(V_\beta,\kappa_\beta,h_\beta)_{\beta\in J}$ the quadratic localization system fixed in Definition~\ref{def:espacios-geometria-acotada-con-frontera}. The construction of $\mathfrak K$ takes place in the same Fermi and interior chart families, through the two fixed linear transformations in the proof of Corollary~\ref{cor:geometria-acotada-frontera-regularidad-uniforme-amann}. Denote the normal rescaling factor by $b_N:=2r_B$. Each support of $\pi_\alpha$ or $h_\beta$ lies in a common shrunken box. The packing argument of Proposition~\ref{prop:atlas-frontera-uniforme}, applied in both directions, gives
\[
 L:=\max\left\{
 \sup_{\alpha\in I}\#\left\{\beta\in J\middle|
 \operatorname{supp}\pi_\alpha\cap\operatorname{supp}h_\beta
 \neq\varnothing\right\},
 \sup_{\beta\in J}\#\left\{\alpha\in I\middle|
 \operatorname{supp}h_\beta\cap\operatorname{supp}\pi_\alpha
 \neq\varnothing\right\}
 \right\}<+\infty.
\]

Fix $k>|s|+n+1$. The $BUC^{k+1}$ bounds for transitions and cutoffs, the upper and lower Jacobian bounds, and Lemma~\ref{lem:operadores-localizados-semiespacio-fermi} provide numbers $M_{\to,\mathcal A},M_{\leftarrow,\mathcal A}>0$, for $\mathcal A=H^{s,p}$ or $F^s_{p,q}$, bounding, respectively, each localized operator passing from $\mathcal T_F$ to Amann's localization system and each operator in the reverse direction. For two boundary charts, the transition preserves the normal parameter within each system. Comparing $\mathcal T_F$ with the rescaled charts in the proof of Corollary~\ref{cor:geometria-acotada-frontera-regularidad-uniforme-amann} introduces the fixed normal factor $b_N$; thus, the transition has the form $(x',t)\mapsto(\vartheta(x'),\lambda t)$, with $\lambda\in\{1,b_N,b_N^{-1}\}$. Use item (b) with $\displaystyle \lambda_-:=\min\{1,b_N,b_N^{-1}\}$ and $\displaystyle \lambda_+:=\displaystyle\max\{1,b_N,b_N^{-1}\}$. On a mixed overlap, the support is at uniform distance from $t=0$ and item (c) is used; between interior charts, use Lemma~\ref{lema: multiplicadores y difeomorfismos sobolev}. In particular, both bounds remain valid for $s<0$. The numbers depend only on $n,s,p,k$, also on $q$ in scale $F$, the Jacobian bounds, the $BUC^{k+1}$ bounds for transitions and cutoffs, and the uniform separation of mixed overlaps.

Insert $\displaystyle \sum_{\beta\in J} h_\beta^2=1$ into every localization $\pi_\alpha u$. If $A_\alpha$ denotes its norm in $H^{s,p}(X_\alpha)$ and $B_\beta$ the corresponding summand of \eqref{eq:norma-H-frontera}, then
\[
 A_\alpha\leq M_{\to,H}
 \sum_{\left\{\beta\in J\,\middle|\,
 \operatorname{supp}\pi_\alpha\cap\operatorname{supp}h_\beta
 \neq\varnothing\right\}}B_\beta.
\]
The inequality $\displaystyle (\displaystyle\sum_{j=1}^La_j)^p\leq L^{p-1}\displaystyle\sum_{j=1}^La_j^p$ and the two overlap counts give
\[
 \|\mathcal S_{\partial,\mathfrak K}u\|_{
 \mathbb H^{s,p}_{\partial}(\mathfrak K)}
 \leq LM_{\to,H}\|u\|_{H^{s,p}(M)}.
\]
In the reverse direction, insert $\displaystyle \sum_{\alpha\in I}\pi_\alpha^2=1$ into each $h_\beta u$ to obtain
\[
 \|u\|_{H^{s,p}(M)}
 \leq LM_{\leftarrow,H}
 \|\mathcal S_{\partial,\mathfrak K}u\|_{
 \mathbb H^{s,p}_{\partial}(\mathfrak K)}.
\]
Consequently, in \eqref{eq:comparacion-amann-frontera-analisis} we may take
\[
 c^\partial_{s,p}:=(LM_{\leftarrow,H})^{-1},
 \qquad
 C^\partial_{s,p}:=LM_{\to,H}.
\]
The same calculation with $F^s_{p,q}(M)$ norms gives
\eqref{eq:comparacion-amann-frontera-analisis-F}
with
\[
 c^F_{s,p,q}:=(LM_{\leftarrow,F})^{-1},
 \qquad
 C^F_{s,p,q}:=LM_{\to,F}.
\]

For synthesis, localizing $\mathcal R_{\partial,\mathfrak K}v$ using $h_\beta$ involves at most $L$ indices $\alpha$. The same three transition cases provide uniform constants $M_{\mathrm{sin},H}>0$ and $M_{\mathrm{sin},F}(q)>0$ such that each summand has norm at most $M_{\mathrm{sin},H}\|v_\alpha\|_{H^{s,p}(X_\alpha)}$ in $H$, and the analogous bound with $M_{\mathrm{sin},F}(q)$ in $F$. Raising to the power $p$, summing first over $\beta$, and counting overlaps in the reverse direction gives
\[
 \|\mathcal R_{\partial,\mathfrak K}v\|_{H^{s,p}(M)}^p
 \leq L^pM_{\mathrm{sin},H}^p
 \sum_{\alpha\in I}\|v_\alpha\|_{H^{s,p}(X_\alpha)}^p.
\]
Thus, we may take $D^\partial_{s,p}:=LM_{\mathrm{sin},H}$ and, in \eqref{eq:comparacion-amann-frontera-sintesis-F}, $D^F_{s,p,q}:=LM_{\mathrm{sin},F}(q)$. The sum is well defined by local finiteness.

The retraction identity is distributional and does not depend on density. Indeed, if $u\in\mathcal D'(M)$, then $\widetilde\chi_\alpha=1$ on a neighborhood of the coordinate support of $\pi_\alpha$ implies, by definition of multiplication and extension of distributions,
\[
 \mathcal E_\alpha\!\left(
 \widetilde\chi_\alpha
 \widetilde{(\pi_\alpha u)\circ\phi_\alpha^{-1}}
 \right)=\pi_\alpha u.
\]
The family is locally finite and $\displaystyle \sum_{\alpha\in I}\pi_\alpha^2=1$; therefore, in $\mathcal D'(M)$,
\[
 \mathcal R_{\partial,\mathfrak K}
 \mathcal S_{\partial,\mathfrak K}u
 =\sum_{\alpha\in I}\pi_\alpha^2u=u.
\]

Finally, this retraction gives global density exactly in the dense ranges. If $u$ belongs to $H^{s,p}(M)$, or to $F^s_{p,q}(M)$ with $q<\infty$, first truncate $\mathcal S_{\partial,\mathfrak K}u$ to finitely many indices and then approximate each component by the restriction of a Schwartz function, as proved using the quotient in Lemma~\ref{lem:operadores-localizados-semiespacio-fermi}. Synthesis of each finite family belongs to $C_c^\infty(M)$ and converges to $u$ by \eqref{eq:comparacion-amann-frontera-sintesis} or \eqref{eq:comparacion-amann-frontera-sintesis-F}. The case $q=\infty$ is outside the density range.
\end{proof}

\begin{remark}[The two localization systems]
Corollary~\ref{cor:geometria-acotada-frontera-regularidad-uniforme-amann} constructs the normalized atlas and its own cutoffs $\pi_\alpha$; it does not identify them with the cutoffs $h_\beta$ of $\mathcal T_F$. Corollary~\ref{cor:puente-funcional-amann-sobolev-frontera} then uses the constants $L$, $M_{\to,\mathcal A}$, and $M_{\leftarrow,\mathcal A}$ to prove that both systems produce the same scales $H$ and $F$ with norms comparable through explicit inequalities. Thus, uniform regularity provides a connection at the function-space level and does not replace the Fermi localization already used in the trace and embedding proofs.
\end{remark}

\begin{proposition}[Exactness under interpolation of the localized scale]
\label{prop:exactitud-interpolacion-F-frontera}
Let $M$ be a Riemannian manifold of bounded geometry with nonempty boundary. Let $1\leq p<\infty$, $a<b$, $0<\rho<1$, and $r=(1-\rho)a+\rho b$. There exist constants $\lambda_{a,b,r},\Lambda_{a,b,r}>0$ such that
\begin{equation}
\label{eq:exactitud-interpolacion-F-frontera}
 \lambda_{a,b,r}\|u\|_{F^r_{p,p}(M)}
 \leq
 \|u\|_{(F^a_{p,p}(M),F^b_{p,p}(M))_{\rho,p}}
 \leq
 \Lambda_{a,b,r}\|u\|_{F^r_{p,p}(M)}.
\end{equation}
The pair $\mathcal S_{\partial,\mathfrak K}$, $\mathcal R_{\partial,\mathfrak K}$ from the preceding corollary is the same for all values of $a,b,r$.
\end{proposition}

\begin{proof}
Write $Y_t:=\mathbb F^t_{p,p,\partial}(\mathfrak K)$. Use the identity $F^t_{p,p}(\mathbb R^n)=B^t_{p,p}(\mathbb R^n)$ and the Besov analysis and synthesis of Lemma~\ref{lem:sintesis-sucesiones-completas}. This synthesis is bounded for every $1\leq p<\infty$, including $p=1$, and its formulas do not depend on $t$. The two operators realize $F^t_{p,p}(\mathbb R^n)$ as a retract of
\[
 b_p^t:=\left\{w=(w_j)_{j\geq0}\in
 \bigl(L^p(\mathbb R^n)\bigr)^{\mathbb N_0}\middle|
 \sum_{j\in\mathbb N_0}2^{jtp}\|w_j\|_{L^p(\mathbb R^n)}^p<+\infty\right\}.
\]
More precisely, there exist operators $\mathfrak a,\mathfrak r$, independent of $t$, such that, for each $t\in\mathbb R$,
\[
 \mathfrak a\colon F^t_{p,p}(\mathbb R^n)\longrightarrow b_p^t,
 \qquad
 \mathfrak r\colon b_p^t\longrightarrow F^t_{p,p}(\mathbb R^n),
 \qquad
 \mathfrak r\mathfrak a=I_{F^t_{p,p}(\mathbb R^n)}.
\]
There also exist $A_t,B_t>0$, depending only on $n,p,t$ and the fixed resolution, such that
\[
 \|\mathfrak a u\|_{b_p^t}\leq A_t\|u\|_{F^t_{p,p}(\mathbb R^n)},
 \qquad
 \|\mathfrak r w\|_{F^t_{p,p}(\mathbb R^n)}\leq B_t\|w\|_{b_p^t}.
\]
In particular, $B_t^{-1}\|u\|_{F^t_{p,p}(\mathbb R^n)}\leq\|\mathfrak a u\|_{b_p^t}
\leq A_t\|u\|_{F^t_{p,p}(\mathbb R^n)}$; no isometric equality between Triebel--Lizorkin and Besov norms is used. In a coordinate $j$,
\[
 K\bigl(\tau,w_j;
 2^{ja}L^p(\mathbb R^n),2^{jb}L^p(\mathbb R^n)\bigr)
 =\min\{2^{ja},\tau2^{jb}\}\|w_j\|_{L^p(\mathbb R^n)}.
\]
If $k(\tau,w)$ is the $\ell^p(I)$ norm of these coordinate functionals, the definitions give
\[
 k(\tau,w)\leq K(\tau,w;b_p^a,b_p^b)\leq2k(\tau,w).
\]
Moreover,
\[
 \int_0^\infty
 \left[\tau^{-\rho}\min\{2^{ja},\tau2^{jb}\}\right]^p
 \frac{d\tau}{\tau}
 =c_{\rho,p}^p2^{jrp},
 \qquad
 c_{\rho,p}^p:=\frac1{p\rho}+\frac1{p(1-\rho)}.
\]
This equality follows by integrating separately before and after $\tau=2^{j(a-b)}$. Tonelli applies because the interpolation index is precisely $p$. Interpolating the same dyadic analysis--synthesis pair gives $m_{\mathbb R},M_{\mathbb R}>0$ such that
\[
 m_{\mathbb R}\|w\|_{F^r_{p,p}(\mathbb R^n)}
 \leq\|w\|_{(F^a_{p,p}(\mathbb R^n),F^b_{p,p}(\mathbb R^n))_{\rho,p}}
 \leq M_{\mathbb R}\|w\|_{F^r_{p,p}(\mathbb R^n)}.
\]

On the half-space, $r_+$ has norm at most one and the same $\mathbf{E}_+$ satisfies $r_+\mathbf{E}_+=I$ at orders $a,b$. If $e_a,e_b$ are its norms, interpolation of $\mathbf{E}_+$ gives the following lower inequality. For the upper inequality, choose $U\in F^r_{p,p}(\mathbb R^n)$ with $r_+U=w$, apply the full-space inequality to $U$, and then take the infimum over all such extensions. Thus, the quotient definition gives
\[
 m_{\mathbb R}e_a^{-(1-\rho)}e_b^{-\rho}
 \|w\|_{F^r_{p,p}(\mathbb R^n_+)}
 \leq
 \|w\|_{(F^a_{p,p}(\mathbb R^n_+),F^b_{p,p}(\mathbb R^n_+))_{\rho,p}}
 \leq M_{\mathbb R}\|w\|_{F^r_{p,p}(\mathbb R^n_+)}.
\]
Taking the outer $\ell^p(I)$ sum and using Tonelli again gives $m_Y,M_Y>0$ with
\begin{equation}
\label{eq:interpolacion-modelo-mixto-frontera}
 m_Y\|v\|_{Y_r}
 \leq\|v\|_{(Y_a,Y_b)_{\rho,p}}
 \leq M_Y\|v\|_{Y_r}.
\end{equation}
One may take $m_Y$ as the minimum of the lower constants for the two models and $M_Y$ as twice the maximum of their upper constants.

Abbreviate $C_t:=C^F_{t,p,p}$ and $D_t:=D^F_{t,p,p}$. Interpolating the common operators, using their retraction identity and \eqref{eq:interpolacion-modelo-mixto-frontera}, gives
\eqref{eq:exactitud-interpolacion-F-frontera}
with
\[
 \lambda_{a,b,r}:=\frac{m_Y}{D_rC_a^{1-\rho}C_b^\rho},
 \qquad
 \Lambda_{a,b,r}:=M_YC_rD_a^{1-\rho}D_b^\rho.
\]
These constants depend only on the indices and uniform atlas data, not on $u$.
\end{proof}

\begin{corollary}[Independence of the Besov endpoints]
\label{cor:independencia-extremos-besov-frontera}
Let $M$ be a Riemannian manifold of bounded geometry with nonempty boundary. Let $1\leq p<\infty$, $1\leq q\leq\infty$, and let $P=(s_0,s_1,\theta)$, $Q=(t_0,t_1,\eta)$ be two triples such that $s_0<s<s_1$, $t_0<s<t_1$, $\theta,\eta\in(0,1)$, and
\[
 s=(1-\theta)s_0+\theta s_1=(1-\eta)t_0+\eta t_1.
\]
There exist constants $c_{P,Q},C_{P,Q}>0$ such that
\[
 c_{P,Q}\|u\|_{B^s_{p,q}(M;s_0,s_1)}
 \leq\|u\|_{B^s_{p,q}(M;t_0,t_1)}
 \leq C_{P,Q}\|u\|_{B^s_{p,q}(M;s_0,s_1)}.
\]
Thus, the space and its topology are independent of the auxiliary endpoints. Henceforth, denote it simply by $B^s_{p,q}(M)$.
\end{corollary}

\begin{proof}
Choose $\displaystyle a<\min\{s_0,t_0\}$ and $\displaystyle b>\displaystyle\max\{s_1,t_1\}$. For $x\in\{s_0,s_1,t_0,t_1\}$, set $\rho_x=\frac{x-a}{b-a}$ and
\[
 Z_x:=(F^a_{p,p}(M),F^b_{p,p}(M))_{\rho_x,p}.
\]
The proposition provides $\lambda_x,\Lambda_x>0$ with
\[
 \lambda_x\|u\|_{F^x_{p,p}(M)}
 \leq\|u\|_{Z_x}
 \leq\Lambda_x\|u\|_{F^x_{p,p}(M)}.
\]
Let $W_s:=(F^a_{p,p}(M),F^b_{p,p}(M))_{\frac{s-a}{b-a},q}$. Theorem~\ref{teo:reiteracion-interpolacion-real} provides $r_P,R_P,r_Q,R_Q>0$ such that
\[
 r_P\|u\|_{W_s}\leq
 \|u\|_{(Z_{s_0},Z_{s_1})_{\theta,q}}
 \leq R_P\|u\|_{W_s},
\]
and the same inequalities with $Q$ in place of $P$. Interpolating the four preceding inequalities gives
\[
 \alpha_P\|u\|_{W_s}
 \leq\|u\|_{B^s_{p,q}(M;s_0,s_1)}
 \leq\beta_P\|u\|_{W_s},
\]
where
\[
 \alpha_P:=\frac{r_P}{\Lambda_{s_0}^{1-\theta}\Lambda_{s_1}^{\theta}},
 \qquad
 \beta_P:=\frac{R_P}{\lambda_{s_0}^{1-\theta}\lambda_{s_1}^{\theta}}.
\]
Define $\alpha_Q,\beta_Q$ by the same formula, with $Q$ in place of $P$. Then we may take
\[
 c_{P,Q}:=\frac{\alpha_Q}{\beta_P},
 \qquad C_{P,Q}:=\frac{\beta_Q}{\alpha_P}.
\]
The outer $\ell^p(I)$ sum was commuted with interpolation only in the intermediate step whose index is $p$; the arbitrary index $q$ enters afterward through reiteration.
\end{proof}

Before taking traces on the half-space, we must verify that the datum does not depend on the chosen full-space representative.

\begin{lemma}[The trace depends only on restriction to the half-space]
\label{lem:traza-desciende-cociente-semiespacio}
Let $1<p<\infty$, $j\in\mathbb N_0$, and $s>j+1/p$. If $U\in H^{s,p}(\mathbb R^n)$ vanishes on $t>0$, then $\operatorname{tr}_1(\partial_t^\ell U)=0$ for $0\leq\ell\leq j$. For $j=0$, the same assertion holds in $F^s_{p,q}(\mathbb R^n)$ and $B^s_{p,q}(\mathbb R^n)$, $1\leq q\leq\infty$, provided that $s>1/p$.
\end{lemma}

\begin{proof}
For $h>0$, set $U_h(x',t)=U(x',t+h)$. This distribution vanishes on $t>-h$. A smooth function $\eta_h(t)$ supported in $(-h/2,h/2)$ and equal to one near zero satisfies $\eta_h\partial_t^\ell U_h=0$. Compatibility with multipliers from Lemma~\ref{lem:compatibilidades-traza-euclidiana} gives $\operatorname{tr}_1(\partial_t^\ell U_h)=0$. Translations converge strongly to the identity on $H^{s,p}(\mathbb R^n)$; continuity of differentiation and the trace permits $h\to 0^{+}$. The same argument applies to $F^s_{p,q}(\mathbb R^n)$ and $B^s_{p,q}(\mathbb R^n)$ when the fine index is finite. For $q=\infty$, choose $0<\varepsilon<s-1/p$ and work in the space of order $s-\varepsilon$ and fine index $1$. Compatibility of traces at different orders identifies the distributional limit with the trace at order $s$.
\end{proof}

\begin{theorem}[Trace and jets on the half-space]
\label{teo:traza-jets-semiespacio}
Let \(1<p<\infty\), \(j\in\mathbb N_0\), and \(s>j+\frac{1}{p}\). The operator
\[
\boldsymbol{\gamma}_jU
:=
\bigl(
U\restriction_{t=0},
\partial_tU\restriction_{t=0},
\dots,
\partial_t^jU\restriction_{t=0}
\bigr)
\]
extends to a retraction
\begin{equation}
\boldsymbol{\gamma}_j\colon
H^{s,p}(\mathbb R^n_+)
\longrightarrow
\prod_{\ell=0}^j
B^{s-\ell-\frac{1}{p}}_{p,p}(\mathbb R^{n-1}).
\label{eq:traza-jets-euclidiana}
\end{equation}
There exists a continuous linear coretraction \(\boldsymbol e_j\) for which \(\boldsymbol{\gamma}_j\boldsymbol e_j=I\). For \(j=0\) and \(1\leq q\leq\infty\), one may replace \(H^{s,p}(\mathbb R^n_+)\) by \(F^s_{p,q}(\mathbb R^n_+)\) and \(B^s_{p,q}(\mathbb R^n_+)\), with target spaces \(B^{s-\frac{1}{p}}_{p,p}(\mathbb R^{n-1})\) and \(B^{s-\frac{1}{p}}_{p,q}(\mathbb R^{n-1})\), respectively.
\end{theorem}

\begin{proof}
Let $u\in H^{s,p}(\mathbb R^n_+)$ and choose an extension $U\in H^{s,p}(\mathbb R^n)$. Define $\gamma_\ell u:=\operatorname{tr}_1(\partial_t^\ell U)$. Two extensions differ by a distribution vanishing on $t>0$; Lemma~\ref{lem:traza-desciende-cociente-semiespacio} shows that the definition is independent of $U$. Continuity of differentiation, the Euclidean trace, and the infimum over extensions give, for $\ell\leq j$,
\[
\|\gamma_\ell u\|_{B^{s-\ell-\frac{1}{p}}_{p,p}(\mathbb R^{n-1})}
\leq C_\ell\inf_{\substack{U\in H^{s,p}(\mathbb R^n)\\r_+U=u}}\|U\|_{H^{s,p}(\mathbb R^n)}
=C_\ell\|u\|_{H^{s,p}(\mathbb R^n_+)}.
\]
Summing the finitely many estimates proves continuity of \(\boldsymbol{\gamma}_j\).

To prescribe the $j+1$ boundary data simultaneously, construct a normal profile for each derivative and estimate the sum of the extensions. Let \((\Delta_k')_{k\geq0}\) be a dyadic resolution in the tangential variables. Choose \(\vartheta\in\mathcal S(\mathbb R)\) with compactly supported Fourier transform and \(\vartheta(0)\neq0\). The finite matrix
\[
\left(
\left.\frac{d^m}{dt^m}\bigl(t^r\vartheta(t)\bigr)\right|_{t=0}
\right)_{0\leq m,r\leq j}
\]
is triangular with diagonal \(r!\vartheta(0)\) and hence invertible. Thus, there exist linear combinations \(\phi_0,\dots,\phi_j\) of the functions \(\vartheta,t\vartheta,\dots,t^j\vartheta\) such that
\begin{equation}
\phi_\ell^{(m)}(0)=\delta_{\ell m},
\qquad 0\leq \ell,m\leq j.
\label{eq:funciones-jets-biortogonales}
\end{equation}
All \(\phi_\ell\) are Schwartz functions with compactly supported Fourier transforms. For a datum \(g\), define on the full space
\begin{equation}
(\mathbf{E}_\ell g)(x',t)
:=
\displaystyle\sum_{k=0}^\infty
2^{-k\ell}\phi_\ell(2^kt)\Delta_k'g(x').
\label{eq:corretraccion-jet-componente}
\end{equation}

We record the discrete estimate controlling superposition in the normal variable. For every \(\phi\in\mathcal S(\mathbb R)\), every nonnegative sequence \((a_k)_{k\geq0}\), and \(1<p<\infty\),
\begin{equation}
\int_{\mathbb R}
\left(
\displaystyle\sum_{k\in\mathbb N_0}2^{\frac{2k}{p}}|\phi(2^kt)|^2a_k^2
\right)^{\frac{p}{2}}dt
\leq C_{\phi,p}\displaystyle\sum_{k\in\mathbb N_0}a_k^p.
\label{eq:hardy-jets-normal}
\end{equation}
Indeed, decompose \(\mathbb R\setminus\{0\}\) into the annuli \(A_m:=\{t\in\mathbb R\mid2^{-m-1}<|t|\leq2^{-m}\}\), \(m\in\mathbb Z\), and extend \(a_k\) by zero for \(k<0\). If \(t\in A_m\), boundedness of \(\phi\) for \(k\leq m\) and its decay of order \(N>\frac{1}{p}\) for \(k>m\) give
\[
2^{\frac{k}{p}}|\phi(2^kt)|
\leq C_N2^{\frac{m}{p}}b_{k-m},
\qquad
b_r:=
\begin{cases}
2^{\frac{r}{p}},&r\leq0,\\
2^{-(N-\frac{1}{p})r},&r>0.
\end{cases}
\]
The sequence \(b\) belongs to \(\ell^1(\mathbb Z)\). Since the \(\ell^2(\mathbb N_0)\) norm does not exceed the \(\ell^1(\mathbb N_0)\) norm, the integral over \(A_m\), whose measure is comparable to \(2^{-m}\), is bounded by \(C(\displaystyle\sum_{k\in\mathbb Z}b_{k-m}a_k)^p\). Summing over \(m\), Young's inequality for sequences from Proposition~\ref{prop:young-sucesiones} proves \eqref{eq:hardy-jets-normal}.

For \(k\geq1\), the Fourier transform of the \(k\)th summand of \eqref{eq:corretraccion-jet-componente} lies in an annulus of radius comparable to \(2^k\): the tangential frequency already lies in that annulus, and the normal frequency lies in a ball of radius \(C2^k\). Finite interaction of blocks and Lemma~\ref{lem:sintesis-bloques-espectrales} give the synthesis estimate
\[
\|\mathbf{E}_\ell g\|_{F^s_{p,2}(\mathbb R^n)}
\leq C
\left\|
\left(
\displaystyle\sum_{k\in\mathbb N_0}
2^{2k(s-\ell)}|\phi_\ell(2^kt)\Delta_k'g(x')|^2
\right)^{\frac{1}{2}}
\right\|_{L^p(\mathbb R^{n-1}\times\mathbb R)}.
\]
Apply \eqref{eq:hardy-jets-normal} for each \(x'\), with \(a_k(x')=2^{k(s-\ell-\frac{1}{p})}|\Delta_k'g(x')|\), and integrate in \(x'\). This gives
\begin{equation}
\|\mathbf{E}_\ell g\|_{H^{s,p}(\mathbb R^n)}^p
\leq C
\displaystyle\sum_{k\in\mathbb N_0}
2^{kp(s-\ell-\frac{1}{p})}\|\Delta_k'g\|_{L^p(\mathbb R^{n-1})}^{p}
=C\|g\|_{B^{s-\ell-\frac{1}{p}}_{p,p}(\mathbb R^{n-1})}^p.
\label{eq:cota-corretraccion-jet-componente}
\end{equation}
In particular, the series in \eqref{eq:corretraccion-jet-componente} converges in \(H^{s,p}(\mathbb R^n)\). Define
\[
\boldsymbol e_j(g_0,\dots,g_j)
:=
\left(\displaystyle\sum_{\ell=0}^j\mathbf{E}_\ell g_\ell\right)\restriction_{\mathbb R^n_+}.
\]
The quotient norm of the restriction space and \eqref{eq:cota-corretraccion-jet-componente} show that this operator is continuous. For Schwartz data, we can differentiate termwise at \eqref{eq:corretraccion-jet-componente}; by \eqref{eq:funciones-jets-biortogonales},
\[
\left.\partial_t^m\mathbf{E}_\ell g\right|_{t=0}
=\delta_{m\ell}\displaystyle\sum_{k\in\mathbb N_0}\Delta_k'g
=\delta_{m\ell}g,
\qquad 0\leq m\leq j.
\]
Density of \(\mathcal S\) in each \(B^{s-\ell-\frac{1}{p}}_{p,p}(\mathbb R^{n-1})\), together with the already proved continuity of the trace operators, extends this identity to all data. Consequently, \(\boldsymbol{\gamma}_j\boldsymbol e_j=I\), and
\[
\|\boldsymbol e_j(g_0,\dots,g_j)\|_{H^{s,p}(\mathbb R^n_+)}
\leq
C\displaystyle\sum_{\ell=0}^j
\|g_\ell\|_{B^{s-\ell-\frac{1}{p}}_{p,p}(\mathbb R^{n-1})}.
\]
For \(j=0\), assertions \(F\) and \(B\) are those of Theorem~\ref{teo:traza-euclidiana-codimension-c} with \(c=1\), restricted to the half-space.
\end{proof}

\begin{theorem}[Trace and normal derivatives at the boundary]
\label{teo:traza-frontera-geometria-acotada}
\index{trace theorem!manifold with boundary}
Let \(M\) be a manifold of dimension \(n\) with nonempty boundary and bounded geometry, and let \(1<p<\infty\), \(j\in\mathbb N_0\), and \(s>j+\frac{1}{p}\). For \(u\in C_c^\infty(M)\), define
\[
\partial_{\boldsymbol{\nu}}^\ell u(x)
:=
\left.\frac{d^\ell}{dt^\ell}
u(\mathcal C(x,t))\right|_{t=0},
\qquad 0\leq\ell\leq j.
\]
Then \[
\boldsymbol{\operatorname{Tr}}_{\partial M,j}u
:=
\bigl(
u\restriction_{\partial M},
\partial_{\boldsymbol{\nu}} u,\dots,\partial_{\boldsymbol{\nu}}^ju
\bigr)
\] extends to a retraction
\begin{equation}
\boldsymbol{\operatorname{Tr}}_{\partial M,j}\colon
H^{s,p}(M)
\longrightarrow
\prod_{\ell=0}^j
B^{s-\ell-\frac{1}{p}}_{p,p}(\partial M).
\label{eq:traza-jets-frontera}
\end{equation}
In particular,
\[
\operatorname{Tr}_{\partial M}\colon
H^{s,p}(M)\longrightarrow B^{s-\frac{1}{p}}_{p,p}(\partial M)
\]
is surjective and has a continuous linear right extension.
\end{theorem}

\begin{proof}
Let \((U_i^\partial,\kappa_i^\partial,h_i)_{i\in I_\partial}\) be the boundary part of the system in Proposition~\ref{prop:atlas-frontera-uniforme}, and write
\[
\kappa_i^{\partial,0}(u):=\kappa_i^\partial(u,0),
\qquad
h_i^\partial:=h_i\restriction_{\partial M}.
\]
By \eqref{eq:cortes-constantes-normal}, on a neighborhood of \(t=0\),
\begin{equation}
\partial_t^\ell\bigl((h_i u)\circ\kappa_i^\partial\bigr)\restriction_{t=0}
=
(h_i^\partial\circ\kappa_i^{\partial,0})
\bigl((\partial_{\boldsymbol{\nu}}^\ell u)\circ\kappa_i^{\partial,0}\bigr).
\label{eq:derivadas-locales-frontera}
\end{equation}
No lower-order Leibniz terms appear; this is the reason for condition \eqref{eq:cortes-constantes-normal}.

Apply \(\boldsymbol{\gamma}_j\) of Theorem~\ref{teo:traza-jets-semiespacio} to each \((h_i u)\circ\kappa_i^\partial\). If \(A(\beta)\) is the set of charts meeting boundary chart \(\beta\), compatibility with multipliers and tangential changes gives, for \(\ell\in\{0,\dots,j\}\),
\[
\bigl\|(h_\beta^\partial\,\partial_{\boldsymbol{\nu}}^\ell u)\circ\kappa_\beta^{\partial,0}\bigr\|_{B^{s-\ell-\frac{1}{p}}_{p,p}(\mathbb R^{n-1})}
\leq C\displaystyle\sum_{i\in A(\beta)}\bigl\|(h_i u)\circ\kappa_i^\partial\bigr\|_{H^{s,p}(\mathbb R^n_+)}.
\]
The cardinality of \(A(\beta)\) is uniformly bounded. Raising to the power \(p\), summing over \(\beta\), and using the fact that each chart appears in uniformly finitely many sets \(A(\beta)\) gives
\begin{equation}
\displaystyle\sum_{\ell=0}^j
\|\partial_{\boldsymbol{\nu}}^\ell u\|_{B^{s-\ell-\frac{1}{p}}_{p,p}(\partial M)}
\leq
C\|u\|_{H^{s,p}(M)}.
\label{eq:cota-jets-frontera}
\end{equation}
This defines the continuous extension of \(\boldsymbol{\operatorname{Tr}}_{\partial M,j}\).

For the coretraction, let \(g=(g_0,\dots,g_j)\) be boundary data and set
\[
g_{\ell,i}:=(h_i^\partial g_\ell)\circ\kappa_i^{\partial,0}.
\]
Choose $\chi_T\in C_c^\infty(B^{n-1}(0,2r_B))$ equal to one on the tangential supports and $\chi_N\in C_c^\infty([0,2r_B))$ equal to one on $[0,r_0]$. Set $\psi(u,t)=\chi_T(u)\chi_N(t)$. Then $\psi\boldsymbol e_j(g_{0,i},\dots,g_{j,i})$ is supported in the chart, and its jet at $t=0$ is the local datum. Transport, extend by zero, and reconstruct:
\begin{equation}
\boldsymbol{\operatorname{Ex}}_{\partial M,j}g
:=
\displaystyle\sum_{i\in I_\partial}
h_i
\left[
\psi\,
\boldsymbol e_j(g_{0,i},\dots,g_{j,i})
\right]\circ(\kappa_i^\partial)^{-1}.
\label{eq:extension-jets-frontera}
\end{equation}
Localizing \(\boldsymbol{\operatorname{Ex}}_{\partial M,j}g\) in chart \(\beta\) involves only the indices in \(A(\beta)\). Uniform bounds for multipliers, chart transitions, and \(\boldsymbol e_j\) bound this localization by the sum of the norms of \(g_{\ell,i}\), with \(i\in A(\beta)\) and \(\ell\in\{0,\dots,j\}\). Raise to the power \(p\) and count overlaps again to obtain
\[
\|\boldsymbol{\operatorname{Ex}}_{\partial M,j}g\|_{H^{s,p}(M)}
\leq
C\displaystyle\sum_{\ell=0}^j
\|g_\ell\|_{B^{s-\ell-\frac{1}{p}}_{p,p}(\partial M)}.
\]
The functions \(h_i\) are constant in the normal direction near the boundary, \(\psi=1\) there, and \(\boldsymbol{\gamma}_j\boldsymbol e_j=I\). Thus, taking the jet of \eqref{eq:extension-jets-frontera} gives \(\displaystyle\sum_{i\in I_\partial}(h_i^\partial)^2g_\ell=g_\ell\) in each component. Hence,
\[
\boldsymbol{\operatorname{Tr}}_{\partial M,j}
\boldsymbol{\operatorname{Ex}}_{\partial M,j}
=I.
\]
\end{proof}

For the noninteger $W$ scale in the following theorem, use the same localized definition as above, with $\mathbb X_i=\mathbb R^n$ in interior charts and $\mathbb X_i=\mathbb R^n_+$ in boundary charts:
\[
 \|u\|_{W^{s,p}_{\mathcal Q}(M)}
 =\left(\sum_{i\in I}\|(h_i u)\circ\kappa_i\|_{W^{s,p}(\mathbb X_i)}^p\right)^{1/p}.
\]
In a boundary chart, extend by zero only across artificial chart faces and
use the restriction norm on the half-space; no extension by zero is made
across the actual boundary. The system independence is proved in
Proposition~\ref{prop:slobodeckij-uniforme-localizacion}.

\begin{theorem}[\(F\), \(B\), and \(W\) traces at the boundary]
\label{teo:traza-BFW-frontera-geometria-acotada}
Let \(M\) be a Riemannian manifold of dimension \(n\) with bounded geometry and nonempty boundary. Let \(1<p<\infty\), \(1\leq q\leq\infty\), and \(s>\frac{1}{p}\). Then
\[
\operatorname{Tr}_{\partial M}\colon
F^s_{p,q}(M)\longrightarrow B^{s-\frac{1}{p}}_{p,p}(\partial M),
\qquad
\operatorname{Tr}_{\partial M}\colon
B^s_{p,q}(M)\longrightarrow B^{s-\frac{1}{p}}_{p,q}(\partial M)
\]
are retractions. If \(s\notin\mathbb N\),
\[
\operatorname{Tr}_{\partial M}\colon
W^{s,p}(M)\longrightarrow B^{s-\frac{1}{p}}_{p,p}(\partial M)
\]
is also a retraction; for \(m\in\mathbb N\), \(m>\frac{1}{p}\), the same holds with \(W^{m,p}(M)\) and \(B^{m-\frac{1}{p}}_{p,p}(\partial M)\).
\end{theorem}

\begin{proof}
Take \(j=0\) in formulas \eqref{eq:traza-jets-frontera} and \eqref{eq:extension-jets-frontera}. If \(u_i=(h_i u)\circ\kappa_i^\partial\), the Euclidean half-space trace, tangential coordinate changes, and uniform overlap count give
\[
\|\operatorname{Tr}_{\partial M}u\|_{B^{s-\frac{1}{p}}_{p,p}(\partial M)}^p
\leq
C\displaystyle\sum_{i\in I_\partial}
\|u_i\|_{F^s_{p,q}(\mathbb R^n_+)}^p
\leq C'\|u\|_{F^s_{p,q}(M)}^p.
\]
For \(g\in B^{s-\frac{1}{p}}_{p,p}(\partial M)\), the same count applied to the reconstruction formula, now with the Euclidean coretraction, gives
\[
\|\operatorname{Ex}_{\partial M,0}g\|_{F^s_{p,q}(M)}
\leq C\|g\|_{B^{s-\frac{1}{p}}_{p,p}(\partial M)}.
\]
Since the cutoffs are constant in the normal direction where the trace is computed, the identity \(\displaystyle\sum_{i\in I_\partial}(h_i^\partial)^2=1\) implies \(\operatorname{Tr}_{\partial M}\operatorname{Ex}_{\partial M,0}=I\).

For the Besov scale, choose \(s_0<s<s_1\) with \(s_0>\frac{1}{p}\) and apply Theorem~\ref{teo:interpolacion-operadores-metodo-K} to the same operators at the endpoints. The interpolation identities
\[
(F^{s_0}_{p,p}(M),F^{s_1}_{p,p}(M))_{\theta,q}=B^s_{p,q}(M)
\]
and \[
(B^{s_0-\frac{1}{p}}_{p,p}(\partial M),
B^{s_1-\frac{1}{p}}_{p,p}(\partial M))_{\theta,q}
=B^{s-\frac{1}{p}}_{p,q}(\partial M)
\] give the second retraction.

If $s\notin\mathbb N$, use restriction $W^{s,p}(\mathbb R^n_+)$ norms in boundary charts, never extension by zero across $t=0$. Proposition~\ref{prop:slobodeckij-uniforme-localizacion} proves that $W^{s,p}(M)=B^s_{p,p}(M)$. This identification reduces the assertion to the Besov retraction with \(q=p\); the target space remains \(B^{s-\frac{1}{p}}_{p,p}(\partial M)\). For $m\in\mathbb N$, the boundary equivalence of Proposition~\ref{prop:independencia-espacios-haz}, applied to the trivial bundle, reduces the assertion to Theorem~\ref{teo:traza-frontera-geometria-acotada} with \(j=0\).
\end{proof}

\section{Version for vector bundles of bounded geometry}

The trace of a section takes values in the restricted bundle. To construct it through components, choose frames adapted to the charts and parallel transport. Bounds on their transition matrices ensure uniform constants in the scalar estimates; compatibility with frame changes allows local traces to be assembled into a section of the restricted bundle.

We use the notion of a bundle of bounded geometry from Definition~\ref{def:haz-vectorial-geometria-acotada}. If \(M\) has boundary, the same expression is understood up to \(\partial M\): require the bound
\eqref{eq:curvatura-haz-acotada}
throughout \(M\), using the induced product connection on \(\Lambda^2T^*M\otimes\operatorname{End}(\mathbf{E})\). Uniform Fermi charts allow these bounds to be formulated with ordinary derivatives up to the boundary.

\begin{proposition}[Uniform synchronous frames]
\label{prop:marcos-sincronos-uniformes}
Let $M$ be a Riemannian manifold of bounded geometry, with or without boundary, and let $\mathbf{E}\to M$ be a bundle of bounded geometry.
\begin{enumerate}[label=(\alph*)]
\item If $M$ has no boundary, its uniform geodesic charts admit synchronous orthonormal frames $\mathbf{e}_i=(\mathbf{e}_{i,1},\dots,\mathbf{e}_{i,r})$ such that the transition matrices, their inverses, and all their derivatives are bounded by constants independent of the indices.
\item If $(M,N)$ is a pair of bounded geometry consisting of manifolds without boundary, the frames can be constructed in the Fermi charts of Proposition~\ref{prop:atlas-fermi-uniforme}: first transport a basis of $\mathbf{E}_{p_i}$ along radial geodesics of $N$, then along normal geodesics of $N$ in $M$. Their restrictions to $N$ are uniform synchronous frames of $\mathbf{E}\restriction_N$.
\item If $M$ has nonempty boundary, in a boundary chart $\kappa_i^\partial(u,t)$ the frame may be chosen parallel along the curves $t\mapsto\kappa_i^\partial(u,t)$. If
$\nabla^{\mathbf{E}}=d+\mathbf{A}_i$
in this frame, then
\begin{equation}
\label{eq:conexion-normal-cero}
\mathbf{A}_i(\boldsymbol{\partial}_t)=0.
\end{equation}
\end{enumerate}
In all three cases, the components of $\mathbf{A}_i$ and all their ordinary derivatives are uniformly bounded.
\end{proposition}

\begin{proof}
Item (a), including the bounds on transition matrices and connection coefficients, is Proposition~\ref{prop:marcos-sincronos-geodesicos-uniformes}. The radial formula producing these bounds was established in the chapter on connections, in Proposition~\ref{prop:formula-curvatura-marco-sincrono}.

For (b), first apply Proposition~\ref{prop:marcos-sincronos-geodesicos-uniformes} to the restricted bundle \(\mathbf{E}\restriction_N\). Its curvature is \[
\mathbf{R}^{\mathbf{E}\restriction_N}(X,Y)=\mathbf{R}^{\mathbf{E}}(X,Y),
\qquad X,Y\in TN,
\], and derivatives for the product connection of \(N\) are finite sums of restrictions of derivatives of \(\mathbf{R}^{\mathbf{E}}\) and contractions with \(\mathbf{II}\). Lemma~\ref{lem:subvariedad-geometria-acotada} and the pair hypotheses therefore provide the required bounds.

From the frame obtained over \(N\), extend each vector by parallel transport along
\[
s\longmapsto
\exp^\perp_{x_i(u)}
\bigl(s\sum_{a=1}^{n-k} z^a\boldsymbol{\nu}_{i,a}(x_i(u))\bigr).
\]
Write frames as columns and use the convention
\[
\nabla^{\mathbf{E}} \mathbf{e}_i=\mathbf{e}_i\boldsymbol{\Omega}_i,
\qquad
\mathscr F_i=d\boldsymbol{\Omega}_i+\boldsymbol{\Omega}_i\wedge\boldsymbol{\Omega}_i.
\]
Transport imposes the normal radial gauge
\[
\sum_{a=1}^{n-k}z^a\Omega_{i,a}(z,u)=0.
\]
Contracting the curvature equation with the radial field gives the exact identities
\begin{align}
\Omega_{i,b}(z,u)
&=
\int_0^1
s\sum_{a=1}^{n-k}z^a\mathscr F_{i,ab}(sz,u)\,ds,
\label{eq:calibre-radial-haz-normal}\\
\Omega_{i,\mu}(z,u)
&=
\Omega_{i,\mu}(0,u)
+
\int_0^1
\sum_{a=1}^{n-k}z^a\mathscr F_{i,a\mu}(sz,u)\,ds.
\label{eq:calibre-radial-haz-tangencial}
\end{align}
Here $k=\dim N$, normal indices $a,b$ belong to $\{1,\ldots,n-k\}$, and $\mu\in\{1,\ldots,k\}$ is tangential. The first identity follows by integrating \[
 \sum_{a=1}^{n-k}z^a\mathscr F_{i,ab}
 =\sum_{a=1}^{n-k}z^a\partial_a\Omega_{i,b}+\Omega_{i,b}
\] along the radius; the second follows from
\[
 \sum_{a=1}^{n-k}z^a\mathscr F_{i,a\mu}
 =\sum_{a=1}^{n-k}z^a\partial_a\Omega_{i,\mu}.
\]
Commutators containing $\displaystyle \sum_{a=1}^{n-k}z^a\Omega_{i,a}$ vanish by radial gauge.

We prove the connection and curvature-component bounds simultaneously. Denote by $A_m$ and $F_m$ common bounds for the supremum norms of partial derivatives of order at most $m$ of $\Omega_{i,A}$ and $\mathscr F_{i,AB}$, respectively, where $i\in I$ and $A,B\in\{1,\ldots,n\}$. At order zero, $F_0$ is controlled by $\|\mathbf R^{\mathbf E}\|_\infty$, uniform comparison of coordinate bases, and orthonormality of the frame. The two integral formulas then give $A_0$, using the frame bounds over $N$.

Suppose $A_{m-1}$ and $F_{m-1}$ have been obtained, with $m\in\mathbb N$. Apply Lemma~\ref{lem:local-expression-higher-order} to the tensor $\mathbf R^{\mathbf E}$ taking values in $\operatorname{End}(\mathbf E)$. Solving for its leading term expresses each partial derivative of order $m$ as the corresponding component of $\nabla^m\mathbf R^{\mathbf E}$ plus a finite sum of partial derivatives of curvature of order less than $m$. The coefficients of this sum use only derivatives of $\Omega_i$ and the Christoffel symbols of order at most $m-1$. The former are controlled by $A_{m-1}$ and the latter by Proposition~\ref{prop:atlas-fermi-uniforme}. This gives $F_m$ from $A_{m-1},F_{m-1}$ and $\|\nabla^m\mathbf R^{\mathbf E}\|_\infty$.

Now differentiate the integral formulas. For $\alpha=(\alpha_z,\alpha_u)\in\mathbb N_0^n$ and $|\alpha|=m$, applying the Leibniz rule to the integrand without the outer factor $s$ gives
\[
 D^\alpha\bigl(z^a\mathscr F_{i,aB}(sz,u)\bigr)
 =z^a s^{|\alpha_z|}(D^\alpha\mathscr F_{i,aB})(sz,u)
  +\alpha_a s^{|\alpha_z|-1}
             (D^{\alpha-e_a}\mathscr F_{i,aB})(sz,u).
\]
Omit the second term if $\alpha_a=0$; thus, no negative power of $s$ is ever integrated. If $r_\perp$ is the common normal radius and $A_m^N$ bounds derivatives of the initial frame over $N$, then
\[
 A_m\leq A_m^N+C_{n,m}\bigl(r_\perp F_m+mF_{m-1}\bigr)+A_{m-1}.
\]
We have obtained $F_m$ before using it to compute $A_m$. The recurrence starts from $A_0,F_0$ and provides bounds at every order, without requiring an order-$m$ connection bound in the step establishing it.

For (c), first choose synchronous orthonormal frames of
$\mathbf E\restriction_{\partial M}$ along boundary geodesics.
Proposition~\ref{prop:frontera-geometria-acotada-inducida} gives bounded
geometry of $\partial M$. The curvature-restriction and covariant-derivative
argument used in (b), now with $N=\partial M$, gives bounded geometry of the
restricted bundle. Part (a) therefore bounds all derivatives of the initial
connection matrices $\Omega_{i,\mu}(u,0)$ uniformly. In a boundary chart,
the same construction along \(t\mapsto\kappa_i^\partial(u,t)\) gives
\[
\boldsymbol{\Omega}_i(\partial_t)=0.
\]
The curvature equation reduces exactly to
\[
\partial_t\Omega_{i,\mu}=\mathscr F_{i,t\mu},
\qquad
\Omega_{i,\mu}(u,t)
=
\Omega_{i,\mu}(u,0)
+\int_0^t\mathscr F_{i,t\mu}(u,\tau)\,d\tau.
\]
For purely tangential derivatives of order $m$, differentiate this integral $m$ times in $u$; the order-$m$ bound for $\mathscr F$ is obtained first by the lemma of Chapter~\ref{cap:sobolev-haces}, exactly as above. If at least one normal derivative occurs, use $\partial_t\Omega_{i,\mu}=\mathscr F_{i,t\mu}$: a total derivative of order $m$ of the connection is expressed through a derivative of order $m-1$ of curvature. The initial case is the order-zero integral and the frame value at $t=0$. These two cases cover every multi-index $\alpha\in\mathbb N_0^n$ and complete the recurrence in boundary charts.

It remains to control the transition matrices. On an overlap, after expressing both frames in the same coordinates, write
\[
\mathbf{e}_i=\mathbf{e}_jG_{ji}.
\]
Since the frames are orthonormal, \(G_{ji}\in O(r)\) in the real case and \(G_{ji}\in U(r)\) in the complex case; thus, \(G_{ji}^{-1}=G_{ji}^{\mathsf T}\) in the former and \(G_{ji}^{-1}=G_{ji}^{*}\) in the latter. If \(\widetilde\Omega_j\) is the pullback of \(\Omega_j\) under the chart transition, compatibility of connections gives the exact equation
\begin{equation}
dG_{ji}
=
G_{ji}\boldsymbol{\Omega}_i-\widetilde\Omega_jG_{ji}.
\label{eq:compatibilidad-transicion-haz-fermi}
\end{equation}
The operator norm of $G_{ji}$ is one. To pass from order $m$ to order $m+1$, let $\alpha\in\mathbb N_0^n$, $|\alpha|=m$, and $A\in\{1,\ldots,n\}$. The Leibniz rule gives
\[
 D^\alpha\partial_A G_{ji}
 =\sum_{\substack{\beta\in\mathbb N_0^n\\\beta\leq\alpha}}
 \binom{\alpha}{\beta}
 \left[
 (D^\beta G_{ji})(D^{\alpha-\beta}\Omega_{i,A})
 -(D^\beta\widetilde\Omega_{j,A})(D^{\alpha-\beta}G_{ji})
 \right].
\]
Every derivative of $G_{ji}$ on the right has order at most $m$. Derivatives of the connection forms are bounded by the preceding step, and those of the pullback by Faà di Bruno and the chart-transition bounds. If $G_m$ is a common bound through order $m$ and $D_m$ controls these forms, we may take
$G_{m+1}=G_m+2^{m+1}D_mG_m$
after incorporating the matrix-norm constants into $D_m$. Starting from $G_0=1$, this bounds every derivative. The identity $G_{ji}^{-1}=G_{ji}^*$, or the transpose in the real case, gives the same bounds for the inverse. At \(z=0\), we recover the synchronous frame constructed over \(N\), also proving the restriction assertion. This concludes the proof.
\end{proof}

\begin{definition}[Function spaces of sections]
\label{def:espacios-secciones-haz-geometria-acotada}
Let $M$ be a Riemannian manifold of bounded geometry, with or without boundary, and let $\mathbf{E}\to M$ be a bundle of bounded geometry. Fix an admissible quadratic localization system \[
 \mathcal Q=(U_i,\phi_i,h_i)_{i\in I}
\] and, on each $U_i$, one of the synchronous frames $\mathbf e_i=(\mathbf e_{i,1},\dots,\mathbf e_{i,r})$ constructed in Proposition~\ref{prop:marcos-sincronos-uniformes}. Let $\mathbf{u}\in\mathcal D'(M,\mathbf{E})$; when there is a boundary, this means $\mathcal D'(\operatorname{Int}(M),\mathbf E|_{\operatorname{Int}(M)})$, as in the scalar definition. In chart $i$, denote by $u_i^a$ the distributional components of $\mathbf{u}$ in frame $\mathbf{e}_i$, defined by Proposition~\ref{prop:representacion-local-distribucion-coordenadas}. Thus, the component of $h_i \mathbf{u}$ is $h_i u_i^a$. Set $\mathbb X_i=\mathbb R^n$ in an interior chart and $\mathbb X_i=\mathbb R^n_+$ in a boundary chart.

For $s\in\mathbb R$ and $1<p<\infty$, define
\[
H^{s,p}(M,\mathbf{E})
:=
\left\{\mathbf{u}\in\mathcal D'(M,\mathbf{E})\middle|\|\mathbf{u}\|_{H^{s,p}(M,\mathbf{E})}<\infty\right\},
\]
where
\begin{equation}
\|\mathbf{u}\|_{H^{s,p}(M,\mathbf{E})}
:=
\left(
\sum_{i\in I}\sum_{a=1}^r
\bigl\|(h_i u_i^a)\circ\phi_i^{-1}\bigr\|_{H^{s,p}(\mathbb X_i)}^p
\right)^{\frac{1}{p}}.
\label{eq:norma-H-haz}
\end{equation}
For $1\leq p<\infty$ and $1\leq q\leq\infty$, define $F^s_{p,q}(M,\mathbf{E})$ analogously using the norm
\[
\left(
\sum_{i\in I}\sum_{a=1}^r
\bigl\|(h_i u_i^a)\circ\phi_i^{-1}\bigr\|_{F^s_{p,q}(\mathbb X_i)}^p
\right)^{\frac{1}{p}}.
\]
Finally, if $s=(1-\theta)s_0+\theta s_1$, with $s_0<s<s_1$, define
\[
B^s_{p,q}(M,\mathbf{E})
:=
\bigl(F^{s_0}_{p,p}(M,\mathbf{E}),F^{s_1}_{p,p}(M,\mathbf{E})\bigr)_{\theta,q}.
\]
On the half-space, always use the restriction quotient norms of Definition~\ref{def:espacios-geometria-acotada-con-frontera}. Definitions on $N$ and $\partial M$ use the restricted frames and constructions without boundary.
\end{definition}

In boundary charts, extensions by zero are taken only across artificial faces. The face representing $\partial M$ belongs to the model $\mathbb R^n_+$ and is treated by the extension operator of Chapter~\ref{cap:regularidad-intermedia}. Fixed coordinate cutoffs leave a common margin from these artificial faces; thus, multiplication and extension constants are the same throughout the family. In an interior chart, apply Lemma~\ref{lem:margen-uniforme-cortes-admisibles}.

\begin{proposition}[Intrinsic character and integer-order equivalence]
\label{prop:independencia-espacios-haz}
Let $M$ be a Riemannian manifold of bounded geometry, with or without boundary, and let $\mathbf{E}\to M$ be a bundle of bounded geometry. Let $\mathcal Q$ and $\widetilde{\mathcal Q}$ be two admissible quadratic localization systems equipped with the synchronous frames of Proposition~\ref{prop:marcos-sincronos-uniformes}. For each scale $A\in\{H,F,B\}$, write $\|\mathbf{u}\|_{A_{\mathcal Q}}$ for the localized norm constructed with $\mathcal Q$. For every admissible set of indices, there exist $c_A,C_A>0$, independent of $\mathbf{u}$, such that
\[
 c_A\|\mathbf{u}\|_{A_{\mathcal Q}}
 \leq\|\mathbf{u}\|_{A_{\widetilde{\mathcal Q}}}
 \leq C_A\|\mathbf{u}\|_{A_{\mathcal Q}}.
\]
The constants depend on the function-space indices, multiplicity, and uniform transition, cutoff, and frame bounds, but not on the charts. For $m\in\mathbb N_0$ and $1<p<\infty$,
$H^{m,p}(M,\mathbf{E})=W^{m,p}(M,\mathbf{E})$
and there exist $c_{m,p},C_{m,p}>0$, independent of $\mathbf{u}$, such that
\[
 c_{m,p}\|\mathbf{u}\|_{W^{m,p}(M,\mathbf{E})}
 \leq\|\mathbf{u}\|_{H^{m,p}(M,\mathbf{E})}
 \leq C_{m,p}\|\mathbf{u}\|_{W^{m,p}(M,\mathbf{E})}.
\]
More precisely, if $N_{m,p}(\mathbf{u})$ denotes the localized norm of \eqref{eq:norma-H-haz}, there exist $c'_{m,p},C'_{m,p}>0$ such that
\[
 c'_{m,p}\|\mathbf{u}\|_{W^{m,p}(M,\mathbf{E})}
 \leq N_{m,p}(\mathbf{u})
 \leq C'_{m,p}\|\mathbf{u}\|_{W^{m,p}(M,\mathbf{E})},
\]
where the covariant norm is
\begin{equation}
\left(\displaystyle\sum_{\ell=0}^m
\|(\nabla^{\mathbf{E}})^\ell \mathbf{u}\|_{L^p(M,T^{(0,\ell)}(TM)\otimes \mathbf{E})}^p
\right)^{\frac{1}{p}}.
\label{eq:norma-covariante-haz}
\end{equation}
Comparison of covariant and partial derivatives holds more generally for $1\leq p<\infty$. If $\mathbf u=\displaystyle\sum_{a=1}^{r_{\mathbf E}}u_i^a\mathbf e_{i,a}$ on $U_i$, set
\begin{equation}
 N^W_{m,p}(\mathbf u)
 :=\left(\sum_{i\in I}\sum_{a=1}^{r_{\mathbf E}}
 \|(h_i u_i^a)\circ\phi_i^{-1}\|_{W^{m,p}(\mathbb X_i)}^p
 \right)^{1/p}.
 \label{eq:norma-local-W-haz-todos-p}
\end{equation}
There exist $c^W_{m,p},C^W_{m,p}>0$ such that
\[
 c^W_{m,p}\|\mathbf u\|_{W^{m,p}(M,\mathbf E)}
 \leq N^W_{m,p}(\mathbf u)
 \leq C^W_{m,p}\|\mathbf u\|_{W^{m,p}(M,\mathbf E)}.
\]
Smooth sections with compact support in $M$, which may meet the boundary, are dense in these spaces for every $1\leq p<\infty$. Identification with $H^{m,p}(M,\mathbf E)$ holds in the range $1<p<\infty$.
\end{proposition}

\begin{proof}
On an overlap, component vectors are related by \(u_i=G_{ij}u_j\). Proposition~\ref{prop:marcos-sincronos-uniformes} bounds all derivatives of \(G_{ij}\), and Lemma~\ref{lema: multiplicadores y difeomorfismos sobolev} controls matrix multiplication and coordinate change simultaneously. Uniform multiplicity and the finite-sum inequality used in Proposition~\ref{prop:independencia-atlas-frontera} produce the constants $c_A,C_A$ in the stated inequalities. The argument holds in \(F\), and real interpolation transfers it to \(B\).

For the integer-order comparison, use Lemma~\ref{lem:local-expression-higher-order} in each regular chart $(U_i,\phi_i)$. If $q\in\mathbb N$, $a\in\{1,\ldots,r_{\mathbf E}\}$, and $(i_1,\ldots,i_q)\in\{1,\ldots,n\}^q$, its formula is
\[
 ((\nabla^{\mathbf E})^q\mathbf u)^a_{i_1\cdots i_q}
 =\partial_{i_1}\cdots\partial_{i_q}u^a
 +\sum_{b=1}^{r_{\mathbf E}}
  \sum_{\substack{\beta\in\mathbb N_0^n\\|\beta|\leq q-1}}
 (B_\beta)^a_{b,i_1\cdots i_q}\partial^\beta u^b.
\]
The coefficients in that lemma use connection derivatives through order $q-1$; their supremum norms are uniform by the frame proposition. Let $P_q$ be the sum of the $L^p(\phi_i(U_i))$ norms of partial derivatives of all components through order $q$, and let $V_q$ be the sum of the $L^p(U_i)$ norms of the pointwise magnitudes of covariant derivatives through that order in the chart considered. Uniform comparison of metrics and density gives $P_0\leq K_0V_0$. If $P_{q-1}\leq K_{q-1}V_{q-1}$, solving for the leading term in the preceding formula gives
\[
 P_q\leq P_{q-1}+C_qV_q+C_qP_{q-1}
 \leq\bigl(C_q+(1+C_q)K_{q-1}\bigr)V_q.
\]
Thus, $K_q=C_q+(1+C_q)K_{q-1}$ can be chosen independently of the chart. The reverse inequality follows from the same formula without solving for the leading term: $V_q\leq C_q'P_q$. Coefficients in both comparisons are those in the lemma of Chapter~\ref{cap:sobolev-haces}, whose bounds are uniform in the trivializations considered.

Apply the comparisons to $h_i\mathbf u$. The Leibniz rule for the product connection controls its derivatives by
\[
 |\nabla^q(h_i\mathbf u)|
 \leq C_q\sum_{b=0}^q
            |\nabla^{q-b}h_i|\,|\nabla^b\mathbf u|,
 \qquad q\in\{0,\ldots,m\}.
\]
The supremum norms of cutoff derivatives and uniform multiplicity allow raising to the power $p$, summing over $i\in I$, and integrating. For the reverse inequality, use $\displaystyle \mathbf u=\displaystyle\sum_{i\in I}h_i(h_i\mathbf u)$ and the same Leibniz rule; uniformly finitely many summands contribute at each point. Equivalence between a sum of finitely many norms and their $\ell^p(I)$ sum adds only factors depending on $m,n,r_{\mathbf E},p$. All these steps hold for $1\leq p<\infty$: only triangle inequalities, the Leibniz rule, and sums of bounded multiplicity were used. They prove both bounds in $N^W_{m,p}$, including $p=1$. For $1<p<\infty$, Theorem~\ref{teo:identificaciones-H-F-W-entero} identifies Euclidean integer-order norms with Bessel norms; on the half-space, apply it after the extension and restriction of Chapter~\ref{cap:regularidad-intermedia}. This provides the two asserted global constants. These triangular formulas also hold for weak derivatives by Lemma~\ref{lema:meyers-serrin-haz-E}. Thus, the comparisons apply directly to each given section without presupposing an approximation already satisfying trace conditions.

We specify the density used later. The series $\displaystyle \mathbf u=\displaystyle\sum_{i\in I} h_i(h_i\mathbf u)$ converges in $W^{m,p}(M,\mathbf E)$: its tail norms are controlled by tails of the localized $\ell^p(I)$ sum. In each of the finitely many charts of a partial sum, extend components to the full space (using $\mathcal E_m$ from Corollary~\ref{teo:extension-sobolev-semiespacio} in a boundary chart), regularize, and multiply by a fixed auxiliary cutoff. Restriction to the half-space gives approximations smooth up to $t=0$. Reconstruction yields sections in $\Gamma_c(M,\mathbf E)$ converging in $W^{m,p}(M,\mathbf E)$. Their supports are not required to avoid the boundary. This approximation holds for every $1\leq p<\infty$ by continuity of translation and regularization in $L^p(\mathbb R^n_+)$, applied to the finitely many derivatives. For $1<p<\infty$, the same construction works in $H^{m,p}(M,\mathbf E)$ and, together with the preceding comparisons, proves equality of the two spaces as subspaces of distributions on $\operatorname{Int}(M)$.
\end{proof}

The integer-order comparison identifies the derivatives entering the norm. At a noninteger order, we still need to compare values in different fibers. The error between parallel transport and component differences is first order in the distance and is therefore integrable in the fractional seminorm.

\begin{proposition}[Uniform Slobodeckij localization]
\label{prop:slobodeckij-uniforme-localizacion}
Let $M$ and $\mathbf E\to M$ have bounded geometry, with or without boundary, let $1<p<\infty$ and $s=m+\sigma$, where $m\in\mathbb N_0$ and $0<\sigma<1$. In a uniform quadratic system with synchronous frames, set
\[
 \|\mathbf u\|_{W^{s,p}_{\mathcal Q}(M,\mathbf E)}
 :=\left(\sum_{i\in I}\sum_{a=1}^{r_{\mathbf E}}
  \bigl\|(h_i u_i^a)\circ\phi_i^{-1}\bigr\|_{W^{s,p}(\mathbb X_i)}^p
       \right)^{1/p}.
\]
In interior charts, extend by zero outside the chart; in boundary charts, extend by zero only across artificial faces and take the restriction norm on the half-space. These norms define a space $W^{s,p}(M,\mathbf E)$ independently of the chosen system, and
\[
 W^{s,p}(M,\mathbf E)=B^s_{p,p}(M,\mathbf E)
 \quad\hbox{with equivalent norms}.
\]
If $M$ has no boundary, this space agrees with the one defined using $\mathbf u\in W^{m,p}(M,\mathbf E)$ and the seminorm
\begin{equation}
 [\nabla_w^m\mathbf u]_{\sigma,p;r}^p
 :=\int_M\int_{0<d_{\mathbf g}(x,y)<r}
 \frac{|\nabla_w^m\mathbf u(x)-P_{y\to x}^{\mathbf E_m}
                    \nabla_w^m\mathbf u(y)|^p}
      {d_{\mathbf g}(x,y)^{n+\sigma p}}
 \,d\lambda_{\mathbf g}(y)\,d\lambda_{\mathbf g}(x),
 \label{eq:slobodeckij-intrinseco-geometria-acotada}
\end{equation}
where $\mathbf E_m=T^{(0,m)}(TM)\otimes\mathbf E$, transport is taken along the unique short geodesic, and $r>0$ is a sufficiently small uniform radius. Replacing $r$ by another positive radius below the uniform radius gives an equivalent norm.
\end{proposition}

\begin{proof}
In each Euclidean model, the identification $W^{s,p}(M,\mathbf E)=B^s_{p,p}(M,\mathbf E)=F^s_{p,p}(M,\mathbf E)$ has constants depending only on $n,s,p$. It also holds for restriction norms on the half-space, by Euclidean extension and the comparisons of Chapter~\ref{cap:regularidad-intermedia}. Proposition~\ref{prop:exactitud-interpolacion-F-frontera} and its version without boundary in the preceding chapter show that the localized space $F^s_{p,p}(\mathbb R^n)$ agrees with $B^s_{p,p}(\mathbb R^n)$. For sections, use the same analysis and synthesis on every component: the rank is finite, matrix transitions are uniform multipliers, and the sums over charts, components, and blocks here have the same exponent $p$. Tonelli and the estimates in those proofs apply without changing their order. This gives equivalence with $B^s_{p,p}(M,\mathbf E)$ and, by Proposition~\ref{prop:independencia-espacios-haz}, independence of the system. This part of the argument includes the boundary case.

Now suppose $\partial M=\varnothing$. Fix normal charts larger than those used by the cutoffs, and choose $r$ smaller than their uniform margin. If one point lies in the support of $h_i$ and its distance from the other is less than $r$, the short geodesic remains inside the larger chart. In the synchronous frame, the transport matrix $G_i(x,y)$ satisfies
\[
 \|G_i(x,y)-I\|+\|G_i(x,y)^{-1}-I\|
 \leq C d_{\mathbf g}(x,y).
\]
Indeed, integrate the transport equation; its coefficients are the bounded connection matrices evaluated on the curve velocity. Coordinate length does not exceed a constant times Riemannian length. Grönwall gives the inequality with $C$ independent of $i$. The same proof applies to each tensor bundle $\mathbf E_m$.

Fix one of these bundles $\mathbf E_m$ and a section $\mathbf v$ of it; write $v_i$ for its component vector. In both directions we obtain
\[
 \begin{aligned}
 |\mathbf v(x)-P_{y\to x}\mathbf v(y)|
 &\leq C\bigl(|v_i(x)-v_i(y)|+d_{\mathbf g}(x,y)|v_i(y)|\bigr),\\
 |v_i(x)-v_i(y)|
 &\leq C\bigl(|\mathbf v(x)-P_{y\to x}\mathbf v(y)|
                   +d_{\mathbf g}(x,y)|\mathbf v(y)|\bigr).
 \end{aligned}
\]
Distances and measures are uniformly comparable to their coordinate counterparts. The contribution of the error terms is bounded by $C\|\mathbf v\|_{L^{p}(M,\mathbf E_m)}^p$, since
\[
 \sup_{y\in M}\int_{d_{\mathbf g}(x,y)<r}
        d_{\mathbf g}(x,y)^{-n+p(1-\sigma)}\,d\lambda_{\mathbf g}(x)
 \leq C\int_0^r t^{p(1-\sigma)-1}\,dt<\infty.
\]
This is the local comparison between the intrinsic seminorm and the component seminorm, now with a uniform constant.

To localize, we use the identity
\[
 h_i(x)\mathbf v(x)-h_i(y)P_{y\to x}\mathbf v(y)
 =h_i(x)\bigl(\mathbf v(x)-P_{y\to x}\mathbf v(y)\bigr)
  +(h_i(x)-h_i(y))P_{y\to x}\mathbf v(y).
\]
We have $|h_i(x)-h_i(y)|\leq C d_{\mathbf g}(x,y)$, and at most $2L$ supports contribute to each pair. We raise to the power $p$, sum, and apply the preceding integral estimate. Conversely, we reconstruct $\displaystyle \mathbf v=\displaystyle\sum_{i\in I} h_i(h_i\mathbf v)$ and use the same finite-sum estimate. The parts of the extension by zero lying outside the chart are separated from the support by a uniform distance $\delta>0$; their contribution is bounded by $C\delta^{-\sigma p}\|\mathbf v\|_{L^{p}(M,\mathbf E_m)}^p$. Thus the local norms and the $\|\mathbf v\|_{L^{p}(M,\mathbf E_m)}+[\mathbf v]_{\sigma,p;r}$ norm are equivalent.

For $s=m+\sigma$, we use the triangular expression for covariant derivatives in Lemma~\ref{lem:local-expression-higher-order}, both directly and to solve for the partial derivatives of order $m$. Each lower-order term contains a derivative of order at most $m-1$ and a uniformly smooth coefficient. It is controlled in $W^{\sigma,p}(\mathbb R^n_+)$ by
$W^{1,p}(\mathbb R^n_+)\hookrightarrow W^{\sigma,p}(\mathbb R^n_+)$
and the multiplier estimate of Lemma~\ref{lem:W1p-controla-slobodeckij-global}. These terms can therefore be absorbed into $\|\mathbf u\|_{W^{m,p}(M,\mathbf E)}$. The comparison at order $\sigma$ applied to $\nabla_w^m\mathbf u$, together with the Leibniz rule for the cutoffs, proves equivalence with $\|\mathbf u\|_{W^{m,p}(M,\mathbf E)}+[\nabla_w^m\mathbf u]_{\sigma,p;r}$. The local formulas hold for weak derivatives; density in the fractional norm has not been assumed.

Finally, if $0<r_1<r_2$ are two admissible radii, the difference between the integrals is taken over $r_1\leq d_{\mathbf g}(x,y)<r_2$. The isometry of parallel transport and the uniform bound on the volume of $B_{\mathbf g}(y,r_2)$ bound this difference by $C_{r_1,r_2}\|\nabla_w^m\mathbf u\|_{L^{p}(M,\mathbf E_m)}^p$. This proves independence of the radius without requiring $M$ to have finite total volume.
\end{proof}

\begin{proposition}[Interpolation and complex Sobolev duality for bundles]
\label{prop:dualidad-sobolev-haces-geometria-acotada}
Let $(M,\mathbf{g})$ be a manifold without boundary of bounded geometry, and let $(\mathbf{E},\mathbf{h}_{\mathbf{E}},\nabla^{\mathbf{E}})\to M$ be a complex Hermitian bundle of bounded geometry; we adopt the convention that $\mathbf{h}_{\mathbf{E}}$ is linear in the first variable. If $s_0<s_1$, $0<\theta<1$, $s_\theta=(1-\theta)s_0+\theta s_1$, and $1<p<\infty$, then
\begin{equation}
\label{eq:interpolacion-compleja-sobolev-haces-geometria-acotada}
 [H^{s_0,p}(M,\mathbf{E}),H^{s_1,p}(M,\mathbf{E})]_\theta
 =H^{s_\theta,p}(M,\mathbf{E})
\end{equation}
with equivalent norms. These identifications are induced by the identity on test sections and are independent of the admissible quadratic localization system and synchronous frames used to define the norms.

If $s\in\mathbb R$, $1<p<\infty$, and $p'=\frac{p}{p-1}$, the complex-bilinear pairing \[
 \mathbf{E}_x\times\overline{\mathbf{E}}_x\longrightarrow\mathbb C,
 \qquad (e,\overline f)\longmapsto \mathbf{h}_{\mathbf{E}}(e,f),
\] followed by integration extends uniquely to a continuous bilinear pairing
\begin{equation}
\label{eq:emparejamiento-sobolev-haces-geometria-acotada}
 \langle\cdot,\cdot\rangle_{s,p;\mathbf{E}}\colon
 H^{s,p}(M,\mathbf{E})\times H^{-s,p'}(M,\overline{\mathbf{E}})
 \longrightarrow\mathbb C.
\end{equation}
The map \begin{equation}
\label{eq:dualidad-sobolev-haces-geometria-acotada}
 \mathfrak J_{s,p,\mathbf{E}}\colon
 H^{-s,p'}(M,\overline{\mathbf{E}})
 \longrightarrow (H^{s,p}(M,\mathbf{E}))',
 \qquad
 \mathfrak J_{s,p,\mathbf{E}}(\mathbf{v})(\mathbf{u}):=\langle \mathbf{u},\mathbf{v}\rangle_{s,p;\mathbf{E}},
\end{equation} is a complex-linear topological isomorphism. Pointwise conjugation also induces a complex-linear isometric isomorphism
\begin{equation}
\label{eq:conjugacion-sobolev-haces-localizados}
 \mathcal C_{s,p,\mathbf{E}}\colon
 \overline{H^{s,p}(M,\mathbf{E})}\longrightarrow H^{s,p}(M,\overline{\mathbf{E}}),
 \qquad \mathcal C_{s,p,\mathbf{E}}(\overline{\mathbf{u}}):=\overline{\mathbf{u}}.
\end{equation}
Thus, equivalently,
\[
 \overline{H^{-s,p'}(M,\mathbf{E})}
 \cong (H^{s,p}(M,\mathbf{E}))'.
\]
In the real case all bars are omitted.

These identifications are compatible with \eqref{eq:interpolacion-compleja-sobolev-haces-geometria-acotada}: when the orders $s_0,s_1$ are interpolated, the spaces in the second entry interpolate between orders $-s_0,-s_1$ and yield order $-s_\theta$ (interchanging the endpoints if necessary); the resulting pairing is the unique extension of the same integral pairing on $\Gamma_c(\mathbf{E})\times\Gamma_c(\overline{\mathbf{E}})$. Pointwise conjugation likewise commutes with these identifications.
\end{proposition}

\begin{proof}
Fix an admissible quadratic localization system $\mathcal Q=(U_i,\phi_i,h_i)_{i\in I}$, with $\displaystyle \sum_{i\in I} h_i^2=1$, and unitary synchronous frames $\mathbf e_i=(\mathbf e_{i,1},\ldots,\mathbf e_{i,r})$. Choose cutoffs $\chi_i\in C_c^\infty(U_i)$ equal to one on a neighborhood of $\operatorname{supp}h_i$, and use the conjugate frame $\overline{\mathbf{e}}_i$ in $\overline{\mathbf{E}}$. If \[
 \mathbf{u}=\sum_{a=1}^r u_i^a \mathbf{e}_{i,a},\qquad
 \mathbf{v}=\sum_{a=1}^r v_i^a\overline{\mathbf{e}}_{i,a},
\] and $\mathbf{v}=\overline{\mathbf{w}}$ is smooth, then $v_i^a=\overline{w_i^a}$ and $\displaystyle \mathbf{h}_{\mathbf{E}}(\mathbf{u},\mathbf{w})=\displaystyle\sum_{a=1}^r u_i^av_i^a$. Set
\[
 a_i(x):=\sqrt{\det\bigl((\phi_i^{-1})^*\mathbf{g}\bigr)(x)}.
\]
The functions $a_i$ and $a_i^{-1}$ and the cutoffs involved satisfy uniform multiplier bounds of every order.

Let $\langle z,\varphi\rangle_{\mathbb R^n}$ denote the Euclidean \emph{bilinear} duality between $H^{-s,p'}(\mathbb R^n)$ and $H^{s,p}(\mathbb R^n)$. For $\mathbf{u}\in H^{s,p}(M,\mathbf{E})$ and $\mathbf{v}\in H^{-s,p'}(M,\overline{\mathbf{E}})$, define
\begin{equation}
\label{eq:formula-local-dualidad-haces-geometria-acotada}
 \langle \mathbf{u},\mathbf{v}\rangle_{s,p;\mathbf{E}}
 :=
 \sum_{i\in I}\sum_{a=1}^r
 \left\langle
   (h_iv_i^a)\circ\phi_i^{-1},
   a_i\,(h_iu_i^a)\circ\phi_i^{-1}
 \right\rangle_{\mathbb R^n}.
\end{equation}
Since multiplication by $a_i$ is uniformly bounded on $H^{s,p}(\mathbb R^n)$, Euclidean duality and Hölder's inequality for sequences give
\[
 |\langle \mathbf{u},\mathbf{v}\rangle_{s,p;\mathbf{E}}|
 \leq C
 \|\mathbf{u}\|_{H^{s,p}(M,\mathbf{E})}
 \|\mathbf{v}\|_{H^{-s,p'}(M,\overline{\mathbf{E}})}.
\]
The series in \eqref{eq:formula-local-dualidad-haces-geometria-acotada} is therefore absolutely convergent. If $\mathbf{u},\mathbf{w}\in\Gamma_c(\mathbf{E})$ and $\mathbf{v}=\overline{\mathbf{w}}$, change of variables and $\displaystyle \sum_{i\in I} h_i^2=1$ give
\[
 \langle \mathbf{u},\overline{\mathbf{w}}\rangle_{s,p;\mathbf{E}}
 =\sum_{i\in I}\int_M h_i^2\mathbf{h}_{\mathbf{E}}(\mathbf{u},\mathbf{w})\,d\lambda_{\mathbf{g}}
 =\int_M\mathbf{h}_{\mathbf{E}}(\mathbf{u},\mathbf{w})\,d\lambda_{\mathbf{g}}.
\]
Smooth compactly supported sections are dense in both spaces: this follows by approximation with finitely supported sequences having Schwartz components in the retraction--coretraction used below. Thus the extension is unique.

Conjugation changes the components $u_i^a$ into $\overline{u_i^a}$ and commutes with the cutoffs and real coordinate changes. Since conjugation is an isometry on each Euclidean space $H^{s,p}(\mathbb R^n,\mathbb C^r)$, the rule \eqref{eq:conjugacion-sobolev-haces-localizados} is a linear isometric isomorphism; its linearity uses $\lambda\overline{\mathbf{u}}=\overline{\overline\lambda \mathbf{u}}$ in the conjugate domain.

We now construct a retraction that will prove interpolation and then the surjectivity of \eqref{eq:dualidad-sobolev-haces-geometria-acotada} without presupposing global duality itself. For each $r_0\in\mathbb R$, write
\[
 \mathbb H_{\mathbf{E}}^{r_0,p}
 :=\ell^p\!\left(I;H^{r_0,p}(\mathbb R^n,\mathbb C^r)\right)
\]
and define, using the same formulas for every $r_0$,
\[
 \mathsf A_{\mathbf{E}}\mathbf{u}:=\bigl((h_i\mathbf{u}_i)\circ\phi_i^{-1}\bigr)_{i\in I},
\]
\[
 \mathsf B_{\mathbf{E}}(Z_i)
 :=\sum_{i\in I}h_i\sum_{a=1}^r
 \left[
   \bigl((\chi_i\circ\phi_i^{-1})Z_i\bigr)\circ\phi_i
 \right]^a\mathbf{e}_{i,a},
\]
where each summand is extended by zero. The uniform multiplier and coordinate-change theorems and finite multiplicity show that \[
 \mathsf A_{\mathbf{E}}\colon H^{r_0,p}(M,\mathbf{E})\longrightarrow\mathbb H_{\mathbf{E}}^{r_0,p},
 \qquad
 \mathsf B_{\mathbf{E}}\colon\mathbb H_{\mathbf{E}}^{r_0,p}\longrightarrow H^{r_0,p}(M,\mathbf{E})
\] are continuous; moreover, $\displaystyle \mathsf B_{\mathbf{E}}\mathsf A_{\mathbf{E}}\mathbf{u}=\displaystyle\sum_{i\in I} h_i^2\mathbf{u}=\mathbf{u}$.

Apply these bounds with $r_0=s_0,s_1$. Theorem~\ref{teo:interpolacion-escalas-BF}, applied to the $r$ components, and Proposition~\ref{prop:interpolacion-ellp-valores-banach} identify
\[
 [\mathbb H_{\mathbf{E}}^{s_0,p},\mathbb H_{\mathbf{E}}^{s_1,p}]_\theta
 =\mathbb H_{\mathbf{E}}^{s_\theta,p}
\]
with equivalent norms. The complex interpolation theorem for operators~\ref{teo:interpolacion-compleja-operadores}, applied to $\mathsf A_{\mathbf{E}}$ and $\mathsf B_{\mathbf{E}}$, preserves the identity $\mathsf B_{\mathbf{E}}\mathsf A_{\mathbf{E}}=I$ and proves \eqref{eq:interpolacion-compleja-sobolev-haces-geometria-acotada}. Since the formulas for both operators are independent of $r_0$, this identification is the identity on test sections. Independence of the quadratic localization system also follows from Proposition~\ref{prop:independencia-espacios-haz}.

Let $\Lambda\in(H^{s,p}(M,\mathbf{E}))'$. Euclidean bilinear duality and $(\ell^p(X_i))'=\ell^{p'}(X_i')$ provide \[
 W=(W_i)_{i\in I}
 \in\ell^{p'}\!\left(I;H^{-s,p'}(\mathbb R^n,\mathbb C^r)\right)
\] such that
\begin{equation}
\label{eq:representante-secuencia-dualidad-haces}
 (\Lambda\circ\mathsf B_{\mathbf{E}})(Z)
 =\sum_{i\in I}\langle W_i,Z_i\rangle_{\mathbb R^n},
 \qquad
 \|W\|_{\ell^{p'}(I;H^{-s,p'}(\mathbb R^n,\mathbb C^r))}
 \leq
 \|\Lambda\|_{(H^{s,p}(M,\mathbf{E}))'}
 \|\mathsf B_{\mathbf{E}}\|_{\mathcal L(\mathbb H_{\mathbf{E}}^{s,p},H^{s,p}(M,\mathbf{E}))}.
\end{equation}
Define an $\overline{\mathbf{E}}$-valued distribution $\mathbf{v}$ by
\begin{equation}
\label{eq:transpuesto-analisis-dualidad-haces}
 \langle \mathbf{u},\mathbf{v}\rangle
 :=\sum_{i\in I}
 \left\langle W_i,(h_i\mathbf{u}_i)\circ\phi_i^{-1}\right\rangle_{\mathbb R^n},
 \qquad \mathbf{u}\in\Gamma_c(\mathbf{E}).
\end{equation}
The sum is finite for each $\mathbf{u}$.

We verify directly, without using the duality being proved, that $\mathbf{v}\in H^{-s,p'}(M,\overline{\mathbf{E}})$. For each $j$, let $\mathsf L_{ij}$ be the local operator that sends the components of a section $\mathbf{u}$ in the frame $\mathbf{e}_j$ to
\[
 \mathsf L_{ij}(\mathbf{u}_j)
 :=\bigl(h_i h_j \mathbf{u}_i\bigr)\circ\phi_i^{-1},
 \qquad U_i\cap U_j\neq\varnothing,
\]
inserting $\chi_j$ before extending by zero. It is a composition of a uniform coordinate change, a uniform transition matrix, and a uniform multiplier. Thus \[
 \mathsf L_{ij}\colon H^{s,p}(\mathbb R^n,\mathbb C^r)
 \longrightarrow H^{s,p}(\mathbb R^n,\mathbb C^r)
\] has norm bounded independently of $i,j$. Its bilinear transpose $\mathsf L_{ij}^{\mathsf T}$ satisfies the same bound between the spaces of order $-s$ and exponent $p'$. If $v_j$ denotes the component vector of $\mathbf{v}$ in the frame $\overline{\mathbf{e}}_j$, definition \eqref{eq:transpuesto-analisis-dualidad-haces} gives, in $\mathcal D'(\mathbb R^n,\mathbb C^r)$,
\begin{equation}
\label{eq:formula-local-transpuesto-analisis-dualidad-haces}
 a_j\,(h_jv_j)\circ\phi_j^{-1}
 =
 \sum_{\substack{i\in I\\U_i\cap U_j\neq\varnothing}}
 \mathsf L_{ij}^{\mathsf T}W_i.
\end{equation}
Indeed, both sides have the same action on every test function: this equality is precisely \eqref{eq:transpuesto-analisis-dualidad-haces} applied to $h_j\mathbf{u}$ and written in the chart $j$. Uniform multiplicity, the bound on the transposes, and multiplication by $a_j^{-1}$ give
\[
 \|\mathbf{v}\|_{H^{-s,p'}(M,\overline{\mathbf{E}})}^{p'}
 \leq C\sum_{j\in I}
 \left(\sum_{\substack{i\in I\\U_i\cap U_j\neq\varnothing}}
 \|W_i\|_{H^{-s,p'}(\mathbb R^n,\mathbb C^r)}\right)^{p'}
 \leq C'\sum_{i\in I}
 \|W_i\|_{H^{-s,p'}(\mathbb R^n,\mathbb C^r)}^{p'}.
\]
Thus $\mathbf{v}$ belongs to the asserted space, and its norm is bounded by a multiple of $\|\Lambda\|_{(H^{s,p}(M,\mathbf{E}))'}$.

Finally, for $\mathbf{u}\in H^{s,p}(M,\mathbf{E})$,
\[
 \Lambda(\mathbf{u})
 =\Lambda(\mathsf B_{\mathbf{E}}\mathsf A_{\mathbf{E}}\mathbf{u})
 =\sum_{i\in I}\langle W_i,(h_i\mathbf{u}_i)\circ\phi_i^{-1}\rangle_{\mathbb R^n}
 =\langle \mathbf{u},\mathbf{v}\rangle_{s,p;\mathbf{E}}.
\]
Therefore $\mathfrak J_{s,p,\mathbf{E}}$ is surjective, and the preceding construction proves continuity of its inverse. If $\mathfrak J_{s,p,\mathbf{E}}(\mathbf{v})=0$, then $\mathbf{v}$ vanishes on $\Gamma_c(\mathbf{E})$; since $\mathbf{h}_{\mathbf{E}}$ identifies $\mathbf{E}$ linearly with $(\overline{\mathbf{E}})^*$, this means that $\mathbf{v}$ is the zero distribution. Thus the map is also injective. Its continuity has already been proved using \eqref{eq:formula-local-dualidad-haces-geometria-acotada}.

Finally, the local formulas defining the pairing, conjugation, $\mathsf A_{\mathbf{E}}$, and $\mathsf B_{\mathbf{E}}$ are the same at every order. Density of $\Gamma_c(\mathbf{E})$ and $\Gamma_c(\overline{\mathbf{E}})$ and uniqueness of the continuous extension identify the pairing at the interpolated order with \eqref{eq:emparejamiento-sobolev-haces-geometria-acotada}. This proves the asserted compatibility with interpolation and conjugation. The representative has been constructed using local duality and the analysis and synthesis operators.
\end{proof}

\begin{theorem}[Identification of heat potentials with the local spaces]
\label{teo:bessel-calor-triebel-local-haces-geometria-acotada}
Let $(M,\mathbf g)$ be a manifold without boundary of bounded geometry, and let $(\mathbf E,\mathbf h_{\mathbf E},\nabla^{\mathbf E})\to M$ be a finite-rank bundle of bounded geometry. For $s\in\mathbb R$ and $1<p<\infty$, the identity on test sections induces the identifications
\begin{equation}
\label{eq:bessel-calor-Fp2-local-haz}
 H_{L_{\mathbf E}}^{s,p}(M,\mathbf E)
 =H^{s,p}(M,\mathbf E)
 =F^s_{p,2}(M,\mathbf E).
\end{equation}
In particular, for any admissible quadratic system $\mathcal Q=(U_i,\phi_i,h_i)_{i\in I}$ with synchronous frames, there exist $c,C>0$ such that
\begin{equation}
\label{eq:normas-calor-local-Fp2-haz}
 c\|\mathbf u\|_{H_{L_{\mathbf E}}^{s,p}(M,\mathbf E)}
 \leq
 \left(\sum_{i\in I}\sum_{a=1}^{r_{\mathbf E}}
 \left\|\widetilde{(h_i u_i^a)\circ\phi_i^{-1}}
                  \right\|_{F^s_{p,2}(\mathbb R^n)}^p\right)^{1/p}
 \leq C\|\mathbf u\|_{H_{L_{\mathbf E}}^{s,p}(M,\mathbf E)}.
\end{equation}
The constants depend on $n,r_{\mathbf E},s,p$, the uniform bounds of the localization system, and finitely many bounded-geometry bounds for $M$ and $\mathbf E$. They do not depend on the chart index, the number of charts, or the section $\mathbf u$. For $s<0$, the first equality also gives a canonical injective realization of the completion defined by heat potentials as a space of $\mathbf E$-valued distributions.
\end{theorem}

\begin{proof}
The Euclidean identification $H^{s,p}(\mathbb R^n)=F^s_{p,2}(\mathbb R^n)$ in Theorem~\ref{teo:identificaciones-H-F-W-entero} has constants depending only on $n,s,p$. Apply it to each component and sum the $p$th powers over $i\in I$ and $a\in\{1,\ldots,r_{\mathbf E}\}$. This proves the second equality in \eqref{eq:bessel-calor-Fp2-local-haz} and its two norm comparisons for every real order.

First consider $s\geq0$. For $s=0$, all three norms are equivalent to the $L^p(M,\mathbf E)$ norm, by uniform comparison of the measures and fiber norms. If $s=m$ is a positive integer, Theorem~\ref{teo:yoshida-bessel-covariante-geometria-acotada} and Proposition~\ref{prop:independencia-espacios-haz} give
\[
 H_{L_{\mathbf E}}^{m,p}(M,\mathbf E)
 =W^{m,p}(M,\mathbf E)=H^{m,p}(M,\mathbf E).
\]
Both equalities preserve the section as a distribution, and the corresponding norms are bounded by one another.

For a positive noninteger order, choose $m\in\mathbb N$ with $m>s$ and set $\theta=s/m\in(0,1)$. At the endpoints $0,m$ we have the same compatible pair of spaces, with equivalent norms. Complex interpolation of the identity and its inverse, by Theorem~\ref{teo:interpolacion-compleja-operadores}, preserves this equivalence. On the one hand, Theorem~\ref{teo:interpolacion-compleja-bessel-laplaciano-completo} gives
\[
 [L^p(M,\mathbf E),H_{L_{\mathbf E}}^{m,p}(M,\mathbf E)]_\theta
 =H_{L_{\mathbf E}}^{s,p}(M,\mathbf E).
\]
On the other hand, Proposition~\ref{prop:dualidad-sobolev-haces-geometria-acotada} identifies the same interpolated space with $H^{s,p}(M,\mathbf E)$. The identity of the initial spaces therefore yields the first equality sought. The imaginary-power bound that makes this step possible is the one in Lemma~\ref{lem:potencias-imaginarias-bochner-haz-completo}, applied to the connection Laplacian of $\mathbf E$.

Negative order remains. Write $s=-a$, with $a>0$, and $p'=p/(p-1)$. The conjugate bundle $\overline{\mathbf E}$ has the same geometric bounds as $\mathbf E$. The positive-order result already proved identifies
\[
 H_{L_{\overline{\mathbf E}}}^{a,p'}(M,\overline{\mathbf E})
 =H^{a,p'}(M,\overline{\mathbf E})
\]
with equivalent norms. Take the complex-linear duals of both sides and use Theorem~\ref{teo:propiedades-escala-bessel-haces} and Proposition~\ref{prop:dualidad-sobolev-haces-geometria-acotada}, respectively. We obtain an isomorphism between $H_{L_{\mathbf E}}^{-a,p}(M,\mathbf E)$ and $H^{-a,p}(M,\mathbf E)$. On test sections this isomorphism is the identity: both dualities agree with the same bilinear integral pairing, using the metric to pair $\mathbf E$ with $\overline{\mathbf E}$. Density of test sections in the localized positive-order space uniquely determines the corresponding distribution. It also proves injectivity of this realization. In the real case conjugations are omitted.

Let us make the uniformity precise. Integer-order equivalence uses bounds on the metric, its inverse, the connection coefficients, and the cutoffs in the larger charts, together with the overlap multiplicity. The internal power estimate for the Bochner chain adds bounds on the potentials $\mathbf V_\ell$ and $\mathbf W_\ell$ of the finite chain corresponding to the integer $m$. These bounds are expressed in terms of the supremum norms of finitely many derivatives of the curvatures. Interpolation adds the constants of the two endpoint norms and those for imaginary powers; the latter depend only on $p$ and the rank. For $s<0$, the same data are used with $a=-s$ and $p'$. None of these steps introduces a sum of constants over all charts: overlaps are always counted using the same uniform bound. This gives the constants in \eqref{eq:normas-calor-local-Fp2-haz}.
\end{proof}

\begin{remark}[Choice of operator in the presence of a boundary]
\label{obs:comparacion-calor-local-haz-frontera}
The hypothesis $\partial M=\varnothing$ in the preceding theorem is necessary to identify the operator without imposing additional conditions. On a manifold with boundary, the localized space $H^{s,p}(M,\mathbf E)$ is defined by restriction in boundary charts. Powers of a Dirichlet or Neumann realization incorporate the conditions of that realization. For example, on the half-space, for the scalar Dirichlet Laplacian on $L^2(\mathbb R^n_+)$,
\[
 \mathcal D(I+L_D)
 =H^{2,2}(\mathbb R^n_+)\cap H^{1,2}_0(\mathbb R^n_+),
\]
whereas the localized norm of order two admits functions in $H^{2,2}(\mathbb R^n_+)$ with nonzero trace. Thus, before comparing a heat scale with the spaces of this chapter in the presence of a boundary, one must fix the realization and describe its trace conditions.
\end{remark}

\begin{theorem}[Trace of sections on a submanifold]
\label{teo:traza-haz-subvariedad-geometria-acotada}
Let \((M,N)\) be a bounded-geometry pair of smooth manifolds without boundary, where \(M\) has dimension \(n\) and \(N\) has dimension \(k\), \(c=n-k\), and let \(\mathbf{E}\to M\) be a vector bundle of bounded geometry. If \(1<p<\infty\) and \(s>\frac{c}{p}\), restriction of smooth sections extends to a retraction
\begin{equation}
\boldsymbol{\operatorname{Tr}}_N^{\mathbf{E}}\colon
H^{s,p}(M,\mathbf{E})\longrightarrow
B^{s-\frac{c}{p}}_{p,p}(N,\mathbf{E}\restriction_N).
\label{eq:traza-haz-subvariedad}
\end{equation}
Moreover, for \(1\leq q\leq\infty\), there are retractions
\begin{equation}
\boldsymbol{\operatorname{Tr}}_N^{\mathbf{E}}\colon F^s_{p,q}(M,\mathbf{E})\longrightarrow B^{s-\frac{c}{p}}_{p,p}(N,\mathbf{E}\restriction_N),
\qquad
\boldsymbol{\operatorname{Tr}}_N^{\mathbf{E}}\colon B^s_{p,q}(M,\mathbf{E})\longrightarrow B^{s-\frac{c}{p}}_{p,q}(N,\mathbf{E}\restriction_N).
\label{eq:trazas-BF-haz-subvariedad}
\end{equation}
Each operator has a continuous linear right inverse.
\end{theorem}

\begin{proof}
In the Fermi frame \(\mathbf{e}_i\), write \(h_i \mathbf{u}=\displaystyle\sum_{a=1}^r u_i^a \mathbf{e}_{i,a}\) and apply \(\operatorname{tr}_c\) to each component. By part (c) of Lemma~\ref{lem:compatibilidades-traza-euclidiana}, the local estimate is the finite sum of the scalar estimates in Theorem~\ref{teo:traza-euclidiana-codimension-c}. On an overlap, parts (a) and (b) of the same lemma give
\begin{equation}
\operatorname{tr}_c(G_{ij}u_j\circ\Phi_{ji})
=
(G_{ij}\restriction_{z=0})
(\operatorname{tr}_c u_j)
\circ(\Phi_{ji}\restriction_{z=0}).
\label{eq:compatibilidad-traza-haz}
\end{equation}
This is the transformation law for the components of a section of \(\mathbf{E}\restriction_N\); hence the local traces define an intrinsic operator. The overlap count in \eqref{eq:cota-traza-global} proves its continuity.

To construct the right inverse, write \(h_i^N \mathbf{f}=\displaystyle\sum_{a=1}^r f_i^a \mathbf{e}_{i,a}\restriction_N\) locally, apply \(\operatorname{ex}_c\) componentwise, and reconstruct by
\begin{equation}
\boldsymbol{\operatorname{Ex}}_{M,N}^{\mathbf{E}}\mathbf{f}
=
\displaystyle\sum_{i\in I_N}\displaystyle\sum_{a=1}^r h_i
\left[\psi\,\operatorname{ex}_c f_i^a\right]
\circ\kappa_i^{-1}\,\mathbf{e}_{i,a}.
\label{eq:extension-haz-subvariedad}
\end{equation}
The uniform bounds on the transition matrices prove continuity as in \eqref{eq:cota-extension-subvariedad}. Upon restriction, \(\operatorname{tr}_c\operatorname{ex}_c=I\) and \(\displaystyle\sum_{i\in I_N}(h_i^N)^2=1\), so $\boldsymbol{\operatorname{Tr}}_N^{\mathbf E}
\boldsymbol{\operatorname{Ex}}_{M,N}^{\mathbf E}=I$. For $F^s_{p,q}(M,\mathbf E)$, use the same formulas with the corresponding Euclidean bounds and sum in $\ell^p(I)$ over the charts. For $q=\infty$, the local formula defines the operator on the entire space, and the compatibility relations are understood in the sense of distributions.

For $B^s_{p,q}(M,\mathbf E)$, choose $s_0<s<s_1$, with $s_0>c/p$, and apply the common operators to the two spaces $F^{s_i}_{p,p}(M,\mathbf E)$. We interpolate these operators, rather than the norms of each chart separately:
\[
 \begin{aligned}
 (F^{s_0}_{p,p}(M,\mathbf E),F^{s_1}_{p,p}(M,\mathbf E))_{\theta,q}
 &=B^s_{p,q}(M,\mathbf E),\\
 (B^{s_0-c/p}_{p,p}(N,\mathbf E|_N),
  B^{s_1-c/p}_{p,p}(N,\mathbf E|_N))_{\theta,q}
 &=B^{s-c/p}_{p,q}(N,\mathbf E|_N).
 \end{aligned}
\]
The second equality follows by reiteration of the $F_{p,p}$ scale. Continuity and the retraction identity are preserved. In particular, $B^s_{p,q}(M,\mathbf E)$ is not identified with an $\ell^p(I)$ sum of local Besov norms when $q\neq p$.
\end{proof}

\begin{theorem}[Trace of sections on a manifold with boundary]
\label{teo:traza-haz-frontera-geometria-acotada}
Let \(M\) be a manifold with nonempty boundary of bounded geometry, and let \(\mathbf{E}\to M\) be a vector bundle of bounded geometry. If \(j\in\mathbb N_0\), \(1<p<\infty\), and \(s>j+\frac{1}{p}\), let \(\widetilde{\boldsymbol{\nu}}\) denote the collar vector field determined by \(\widetilde{\boldsymbol{\nu}}_{\mathcal C(x,t)}=\boldsymbol{\partial}_t\mathcal C(x,t)\). Then \begin{equation}
\boldsymbol{\operatorname{Tr}}^{\mathbf{E}}_{\partial M,j}\mathbf{u}:=\bigl(\mathbf{u}\restriction_{\partial M},
\nabla_{\widetilde{\boldsymbol{\nu}}}^{\mathbf{E}}\mathbf{u}\restriction_{\partial M},
\dots,
(\nabla_{\widetilde{\boldsymbol{\nu}}}^{\mathbf{E}})^j\mathbf{u}\restriction_{\partial M}
\bigr)
\label{eq:jet-normal-haz}
\end{equation} extends to a retraction
\begin{equation}
\boldsymbol{\operatorname{Tr}}^{\mathbf{E}}_{\partial M,j}\colon H^{s,p}(M,\mathbf{E})\longrightarrow
\prod_{\ell=0}^j
B^{s-\ell-\frac{1}{p}}_{p,p}
(\partial M,\mathbf{E}\restriction_{\partial M}).
\label{eq:traza-haz-frontera}
\end{equation}
For \(j=0\) and \(1\leq q\leq\infty\), one also obtains the retractions
\begin{equation}
\boldsymbol{\operatorname{Tr}}^{\mathbf{E}}_{\partial M,0}\colon F^s_{p,q}(M,\mathbf{E})\longrightarrow
B^{s-\frac{1}{p}}_{p,p}(\partial M,\mathbf{E}\restriction_{\partial M}),
\qquad
\boldsymbol{\operatorname{Tr}}^{\mathbf{E}}_{\partial M,0}\colon B^s_{p,q}(M,\mathbf{E})\longrightarrow
B^{s-\frac{1}{p}}_{p,q}(\partial M,\mathbf{E}\restriction_{\partial M}).
\label{eq:trazas-BF-haz-frontera}
\end{equation}
All these operators have continuous linear right inverses.
\end{theorem}

\begin{proof}
Choose the frames in part (c) of Proposition~\ref{prop:marcos-sincronos-uniformes}. By \eqref{eq:conexion-normal-cero}, if \(\mathbf{u}=\displaystyle\sum_{a=1}^r u_i^a\mathbf{e}_{i,a}\), then
\begin{equation}
(\nabla_{\widetilde{\boldsymbol{\nu}}}^{\mathbf{E}})^\ell \mathbf{u}
=\displaystyle\sum_{a=1}^r(\partial_t^\ell u_i^a)\mathbf{e}_{i,a}
\quad\text{in the boundary chart}.
\label{eq:jet-haz-componente}
\end{equation}
On an overlap of two such frames, the connection transformation formula gives $\partial_tG_{ij}=0$, since both normal connection matrices vanish. Consequently, the jet components transform by $G_{ij}|_{t=0}$ without lower-order terms. Moreover, the cutoffs $h_i$ are constant in $t$ near the boundary by \eqref{eq:cortes-constantes-normal}; hence no Leibniz terms appear when localizing the jets. Applying the retraction \(\boldsymbol{\gamma}_j\) of Theorem~\ref{teo:traza-jets-semiespacio} componentwise, summing over the overlaps, and using the bounds on the matrices \(G_{ij}\restriction_{\partial M}\) yields
\[
\displaystyle\sum_{\ell=0}^j
\|(\nabla_{\widetilde{\boldsymbol{\nu}}}^{\mathbf{E}})^\ell \mathbf{u}\restriction_{\partial M}\|_{B^{s-\ell-\frac{1}{p}}_{p,p}(\partial M,\mathbf{E}\restriction_{\partial M})}
\leq C\|\mathbf{u}\|_{H^{s,p}(M,\mathbf{E})}.
\]

The coretraction is obtained by applying \(\boldsymbol e_j\) to the component vectors of the data, multiplying by the fixed cutoff \(\psi\), and summing with factors \(h_i\) using formula \eqref{eq:extension-jets-frontera}. The cutoff \(\psi\) equals one in the region where the jets are computed and allows each summand to be extended by zero outside its chart. The identities \eqref{eq:jet-haz-componente}, \(\boldsymbol{\gamma}_j\boldsymbol e_j=I\), and \(\displaystyle\sum_{i\in I_\partial}(h_i^\partial)^2=1\) show that the jet of the extension is the prescribed datum. The continuity estimates are the scalar estimates summed over finitely many components and a cover of uniform multiplicity. For $j=0$, the estimate in $F^s_{p,q}(M,\mathbf E)$ follows directly from the Euclidean trace and the same counting argument in $\ell^p(I)$. The formula remains valid for $q=\infty$ by the distributional interpretation of the trace. For Besov spaces, fix $s_0<s<s_1$ with $s_0>1/p$ and interpolate the same operators between the global spaces $F^{s_i}_{p,p}(M,\mathbf E)$ and $B^{s_i-1/p}_{p,p}(\partial M,\mathbf E|_{\partial M})$. The definition of Besov spaces and reiteration, exactly as in the preceding theorem, give \eqref{eq:trazas-BF-haz-frontera} without interchanging an $\ell^p(I)$ sum over charts with the fine index $q$. Independence of the frames follows from Proposition~\ref{prop:independencia-espacios-haz} and \eqref{eq:compatibilidad-traza-haz}.
\end{proof}

\begin{remark}[Covariant traces and normal jets]
\label{obs:relacion-trazas-covariantes-jets-normales}
The full covariant trace \(\boldsymbol{\operatorname{Tr}}_\ell \mathbf{u}
=((\nabla^{\mathbf{E}})^\ell \mathbf{u})\restriction_{\partial M}\) takes values in \(\bigl(T^{(0,\ell)}(TM)\otimes \mathbf{E}\bigr)\restriction_{\partial M}\). Write $\widetilde{\boldsymbol{\nu}}$ for the unit normal field extended along the collar geodesics. Since $\nabla^M_{\widetilde{\boldsymbol{\nu}}}\widetilde{\boldsymbol{\nu}}=0$, the tensor differentiation rule of Chapter~\ref{cap:sobolev-haces} gives, for $\ell\in\mathbb N_0$,
\[
 (\nabla^{\ell+1}\mathbf u)(\widetilde{\boldsymbol{\nu}}^{\otimes(\ell+1)})
 =\nabla_{\widetilde{\boldsymbol{\nu}}}^{\mathbf E}
       \bigl((\nabla^\ell\mathbf u)(\widetilde{\boldsymbol{\nu}}^{\otimes\ell})\bigr)
 -\sum_{a=1}^{\ell}(\nabla^\ell\mathbf u)
 (\widetilde{\boldsymbol{\nu}},\ldots,
  \underbrace{\nabla^M_{\widetilde{\boldsymbol{\nu}}}\widetilde{\boldsymbol{\nu}}}_{\text{position }a},
  \ldots,\widetilde{\boldsymbol{\nu}}).
\]
The sum vanishes. The case $\ell=0$ is the definition of $\nabla_{\widetilde{\boldsymbol{\nu}}}^{\mathbf E}\mathbf u$; if the contraction of order $\ell$ agrees with the iterated normal derivative, the preceding equality proves the same at order $\ell+1$. Thus
\[
 (\nabla^\ell\mathbf u)(\widetilde{\boldsymbol{\nu}}^{\otimes\ell})
 =(\nabla_{\widetilde{\boldsymbol{\nu}}}^{\mathbf E})^\ell\mathbf u,
 \qquad \ell\in\mathbb N_0.
\]
The normal jet is therefore the contraction of the full trace with \(\boldsymbol{\nu}^{\otimes\ell}\).

Let us make the converse assertion precise. If \(\mathbf{S}\) is an \(\mathbf{E}\)-valued covariant tensor field, the induced connections satisfy
\begin{equation}
\begin{split}
&\bigl[(\nabla_{\mathbf X}\nabla_{\mathbf Y}
 -\nabla_{\mathbf Y}\nabla_{\mathbf X}
 -\nabla_{[\mathbf X,\mathbf Y]})\mathbf S\bigr]
 (\mathbf V_1,\ldots,\mathbf V_r)\\
&\quad=\mathbf R^{\mathbf E}(\mathbf X,\mathbf Y)
       \bigl(\mathbf S(\mathbf V_1,\ldots,\mathbf V_r)\bigr)\\
&\qquad-\sum_{a=1}^r\mathbf S
 (\mathbf V_1,\ldots,\mathbf R^M(\mathbf X,\mathbf Y)\mathbf V_a,
                      \ldots,\mathbf V_r).
\end{split}
\label{eq:conmutacion-tensor-haz-frontera}
\end{equation}
This identity follows directly from the definitions of the induced connections and from \(\mathbf{R}^M,\mathbf{R}^{\mathbf{E}}\). Moreover, for \(\mathbf{X}\) tangent to the boundary, the Gauss and Weingarten formulas in Proposition~\ref{prop:gauss-weingarten-subvariedad} give
\[
\nabla_{\mathbf{X}}^M\mathbf{Y}=\nabla_{\mathbf{X}}^{\partial M}\mathbf{Y}+\mathbf{II}(\mathbf{X},\mathbf{Y}),
\qquad
\nabla_{\mathbf{X}}^M\boldsymbol{\nu}=-\mathbf{A}_{\boldsymbol{\nu}}\mathbf{X}.
\]

Now suppose that the normal jets of orders $\ell\in\{0,\ldots,m-1\}$ vanish. At order zero this says that $\mathbf u\restriction_{\partial M}=0$. Let $\ell\in\{1,\ldots,m-1\}$, and suppose that all covariant traces of orders less than $\ell$ vanish. The component of $\nabla^\ell\mathbf u$ with all arguments normal is zero by the preceding identity. In any other component, choose a tangential argument $\mathbf X$ and move it to the first position by finitely many adjacent transpositions.

For each transposition, apply \eqref{eq:conmutacion-tensor-haz-frontera} to the corresponding lower-order tensor. The outer derivatives are distributed over its products by Proposition~\ref{prop: leibniz operador *}. The added terms are contractions of the form
\[
 \nabla^a\mathbf R^M*\nabla^{\ell-2-a}\mathbf u
 \quad\text{or}\quad
 \nabla^a\mathbf R^{\mathbf E}*\nabla^{\ell-2-a}\mathbf u,
 \qquad a\in\{0,\ldots,\ell-2\}.
\]
For $\ell=1$ there are no transpositions or curvature terms. For $\ell\geq2$, all factors containing $\mathbf u$ have order less than $\ell$, and their traces vanish by the induction hypothesis. In the remaining component, the tangential derivative is in the first position. Evaluating the derivative of this tensor gives the tangential derivative of a trace of order $\ell-1$ and terms obtained by differentiating its arguments. The former is zero; the latter, described by Gauss and Weingarten, are components of the same trace of order $\ell-1$ and also vanish. Thus all components vanish at order $\ell$, completing the induction up to $m-1$. The reverse implication follows by contracting each full trace with $\boldsymbol\nu^{\otimes\ell}$.

To pass to $\mathbf u\in W^{m,p}(M,\mathbf E)$, with $1\leq p<\infty$, first observe that Theorem~\ref{teo:traza-semiespacio}, applied to the components of $\nabla^\ell\mathbf u$ in the boundary models, and comparison \eqref{eq:norma-local-W-haz-todos-p} give, on setting $\mathbf E_\ell:=T^{(0,\ell)}(TM)\otimes\mathbf E$,
\begin{equation}
 \|\boldsymbol{\operatorname{Tr}}_\ell\mathbf u\|_
 {L^p(\partial M,\mathbf E_\ell|_{\partial M})}
 \leq C_{m,p}\|\mathbf u\|_{W^{m,p}(M,\mathbf E)},
 \qquad 0\leq\ell<m.
 \label{eq:traza-covariante-Lp-geometria-acotada}
\end{equation}
Indeed, each tensor $\nabla^\ell\mathbf u$ belongs to $W^{1,p}(M,T^{(0,\ell)}(TM)\otimes\mathbf E)$, the charts and cutoffs satisfy uniform bounds, and the sum of the $p$th powers of the local estimates is controlled by the multiplicity of the cover. The constant depends on $m,p$, the dimension and rank, and the uniform bounds on the geometry, frames, and cutoffs up to the orders used; it does not depend on the section. This argument includes $p=1$ and, in dimension one, counting measure on the boundary. Contraction with $\boldsymbol\nu^{\otimes\ell}$ yields continuous normal jets taking values in $L^p(\partial M,\mathbf E|_{\partial M})$. For $p>1$, they agree with the Besov jets already constructed, since they agree on smooth sections and both operators are continuous.

To justify the identities in the Sobolev case, fix an open subset
$V\subseteq\partial M$ whose closure is compact and is contained in
a boundary chart. Choose a cutoff $\chi\in C_c^\infty(M)$ supported
in that chart and equal to one on a neighborhood of $\overline V$
in $M$.
Theorem~\ref{teo:meyers-serrin-haz-frontera} provides sections
$\mathbf u_k$ smooth up to the boundary that converge to
$\chi\mathbf u$ in $W^{m,p}(M,\mathbf E)$.
Since the cutoff is constant near $V$, continuity in
\eqref{eq:traza-covariante-Lp-geometria-acotada} implies, for
$0\leq\ell<m$,
\[
 T_{k,\ell}:=
   (\boldsymbol{\operatorname{Tr}}_\ell\mathbf u_k)|_V
 \longrightarrow
 T_\ell:=(\boldsymbol{\operatorname{Tr}}_\ell\mathbf u)|_V
 \quad\text{in }L^p(V,\mathbf E_\ell|_V).
\]
The jets of the approximations may be nonzero.
Consider the commutation and tangential differentiation identities
for smooth sections before imposing the vanishing of the jets.
Each term in those identities is a tangential derivative of some
$T_{k,\ell}$, or a product of such a term with smooth coefficients.
The curvature terms also have this form: the curvature derivatives
are fixed coefficients and the derivative of $\mathbf u_k$
appearing there has order less than $m$.

We make the passage to the limit precise. Let $D$ be a tangential
differential operator with smooth coefficients, from
$\mathbf E_\ell|_V$ to one of the tensor bundles involved, and let
$\boldsymbol\varphi$ be a smooth compactly supported test section
on $V$ of its target bundle. Denote its formal adjoint by $D^*$
with respect to the induced bundle metrics and
$d\lambda_{\widetilde{\mathbf g}}$.
The definition of distributional differentiation and Hölder's
inequality give
\[
 \begin{aligned}
 |\langle D(T_{k,\ell}-T_\ell),\boldsymbol\varphi\rangle|
 &=\left|\int_V
       \langle T_{k,\ell}-T_\ell,D^*\boldsymbol\varphi
       \rangle_{\mathbf h_\ell}\,
       d\lambda_{\widetilde{\mathbf g}}\right|\\
 &\leq
   \|T_{k,\ell}-T_\ell\|_{L^p(V,\mathbf E_\ell|_V)}
   \|D^*\boldsymbol\varphi\|_{L^{p'}(V,\mathbf E_\ell|_V)}
   \longrightarrow0,
 \end{aligned}
\]
where $\mathbf h_\ell$ is the induced metric on $\mathbf E_\ell$,
$p'$ is the conjugate exponent, and $p'=\infty$ if $p=1$.
The second factor is finite because $D^*\boldsymbol\varphi$ is
smooth and compactly supported.
Multiplication by a smooth coefficient can be incorporated into
$D$; hence every term in the identities converges in distributions
on $V$.

Now apply these identities to the limits $T_\ell$.
If the normal jets of orders less than $m$ vanish for $\mathbf u$,
the induction described for smooth sections gives $T_\ell=0$
as a distribution for every $\ell<m$.
Since $T_\ell$ is represented by an $L^p$ section, uniqueness of
the representation of a distribution by a locally integrable
function implies $T_\ell=0$ almost everywhere.
The reverse implication follows by contracting each trace with
$\boldsymbol\nu^{\otimes\ell}$; this contraction is continuous in
$L^p$ and agrees with the normal jet on smooth sections, hence
also on their Sobolev limits.
The sets $V$ cover the boundary, proving the global equivalence
between vanishing of the normal jets and vanishing of all
covariant traces of orders less than $m$.
\end{remark}

This equivalence describes the boundary conditions. To identify $W^{m,p}_0(M,\mathbf E)$, it remains to check that there is no obstruction to approximation at infinity either. Summability of the localizations and the uniform cutoffs also reduce this question to the half-space.

\begin{theorem}[Characterization of $W^{m,p}_0(M,\mathbf E)$ in bounded geometry]
\label{teo:W0-geometria-acotada-jets}
Let $M$ be a manifold of bounded geometry with nonempty boundary, and let $\mathbf{E}\to M$ be a vector bundle of bounded geometry equipped with a bundle metric $\mathbf{h}_{\mathbf{E}}$ and a compatible connection $\nabla^{\mathbf{E}}$; in the complex case the metric is Hermitian. If $m\in\mathbb N$ and $1\leq p<\infty$, then
\[
 W^{m,p}_0(M,\mathbf E)
 =\ker\boldsymbol{\operatorname{Tr}}^{\mathbf E}_{\partial M,m-1}
 =\left\{\mathbf u\in W^{m,p}(M,\mathbf E)\middle|
    \boldsymbol{\operatorname{Tr}}_\ell\mathbf u=0,
                  \ 0\leq\ell<m\right\}.
\]
Here $W^{m,p}_0(M,\mathbf E)$ is the closure, in the $W^{m,p}(M,\mathbf E)$ norm, of smooth sections whose support is a compact subset of $\operatorname{Int}(M)$. For $p=1$, traces and jets are interpreted using the $L^1(\partial M,\mathbf E\restriction_{\partial M})$-valued operators of the preceding remark; no identification with the Bessel scale is required.
\end{theorem}

\begin{proof}
Smooth sections compactly supported in the interior have all traces equal to zero. Continuity of the jet operator established in \eqref{eq:traza-covariante-Lp-geometria-acotada} proves the inclusion of $W^{m,p}_0(M,\mathbf E)$ in the kernel. The last equality in the statement is the preceding remark. It therefore suffices to show that a section $\mathbf u$ whose normal jets of all orders less than $m$ vanish belongs to $W^{m,p}_0(M,\mathbf E)$.

Take the quadratic partition $(h_i)$ and the normally parallel frames used in the trace theorem. Write $\displaystyle h_i\mathbf u=\displaystyle\sum_{a=1}^{r_{\mathbf E}} v_i^a\mathbf e_{i,a}$ and extend each $v_i^a\circ\phi_i^{-1}$ by zero across the artificial faces of the chart. In a boundary chart this yields a function in $W^{m,p}(\mathbb R^n_+)$ with bounded support. Since the cutoff is constant in the normal variable near the boundary and the frame is parallel in that direction, its traces satisfy
\[
 \gamma_0\partial_t^\ell(v_i^a\circ\phi_i^{-1})=0,
 \qquad 0\leq\ell<m.
\]
The trace commutes with tangential derivatives: the identity is proved for smooth functions and passes to the limit in distributions by continuity of the traces. Consequently, we also have $\gamma_0D^\beta(v_i^a\circ\phi_i^{-1})=0$ for every $|\beta|\leq m-1$.

The half-space criterion in Theorem~\ref{teo:traza-semiespacio} implies that each component belongs to $W^{m,p}_0(\mathbb R^n_+)$. We can approximate it in $W^{m,p}(\mathbb R^n_+)$ by functions in $C_c^\infty(\mathbb R^n_+)$. Multiply these approximations by a fixed cutoff equal to one on the component's support and supported inside the artificial faces of the chart. This multiplication preserves convergence and allows the functions to be reconstructed on $M$. Each reconstructed approximation is compactly supported in the interior. In interior charts, proceed in the same way using Euclidean regularization. Local equivalence of norms and finite rank thus show that $h_i\mathbf u\in W^{m,p}_0(M,\mathbf E)$ for each $i$.

To combine the approximations, exhaust the countable index set by finite sets $I_J$ and define
\[
 \mathbf u_J:=\sum_{i\in I_J}h_i(h_i\mathbf u).
\]
Each summand belongs to $W^{m,p}_0(M,\mathbf E)$, since multiplication by $h_i$ is continuous and preserves the interior support of the approximations. By the synthesis estimate and integer-order equivalence of norms,
\[
 \|\mathbf u-\mathbf u_J\|_{W^{m,p}(M,\mathbf E)}
 \leq C\left(\sum_{i\in I\setminus I_J}\sum_{a=1}^{r_{\mathbf E}}
       \|v_i^a\circ\phi_i^{-1}\|_{W^{m,p}(\mathbb X_i)}^p
          \right)^{1/p}\longrightarrow0.
\]
The right-hand side is the tail of a convergent series. Since $W^{m,p}_0(M,\mathbf E)$ is closed, it follows that $\mathbf u\in W^{m,p}_0(M,\mathbf E)$. This proof does not use compactness of $M$.
\end{proof}
\section{Sobolev embeddings in bounded geometry}
\label{sec:encajes-sobolev-geometria-acotada}
\index{Sobolev embedding theorem@Sobolev embedding theorem!in bounded geometry}

Without boundary, Corollary~\ref{cor:puente-funcional-amann-sobolev} formulates the local transfer using a uniformly regular atlas, with constants $c_{s,p},C_{s,p},D_{s,p}$. With boundary, Corollary~\ref{cor:puente-funcional-amann-sobolev-frontera} gives the analysis--synthesis pair on the mixed models $\mathbb R^n$ and $\mathbb R^n_+$, with constants $c^\partial_{s,p},C^\partial_{s,p},D^\partial_{s,p}$. The arguments in this section retain the admissible quadratic localization system fixed when defining the global spaces---with geodesic charts when there is no boundary, and geodesic and Fermi charts when the boundary is nonempty---because it simultaneously incorporates the bundle frames and control of the boundary face. Amann's language thus provides a functional bridge without replacing the covariant estimates already established.

On a compact manifold, local Sobolev estimates are combined using a finite cover. On a noncompact manifold, finiteness must be replaced by three uniform controls: the charts have a common model and radius, the constants of the local operators are independent of the chart, and the cover has bounded multiplicity. The outer sum defining the global spaces also retains the spatial exponent. Thus, when the integrability exponent increases from $p_0$ to $p_1$, the inclusion $\ell^{p_0}(I)\hookrightarrow\ell^{p_1}(I)$ appears and $p_0\leq p_1$ is required. This restriction does not occur on a compact manifold, which has finite measure.

\begin{lemma}[Uniform transfer of a local embedding]
\label{lem:transferencia-encaje-local-geometria-acotada}
Let $M$ be a Riemannian manifold of bounded geometry, with or without boundary, and let $\mathbf{E}\to M$ be a bundle of bounded geometry. Let $1<p_0\leq p_1<\infty$ and, for $\nu\in\{0,1\}$, let $A_\nu(\mathbb R^n)$ be one of the spaces \[
H^{s_\nu,p_\nu}(\mathbb R^n)
\qquad\text{or}\qquad
F^{s_\nu}_{p_\nu,q_\nu}(\mathbb R^n),
\], with $1\leq q_\nu\leq\infty$ in the second case. Suppose that the identity inclusion satisfies
\begin{equation}
\label{eq:encaje-local-modelo-geometria-acotada}
\|v\|_{A_1(\mathbb R^n)}
\leq
C_{\mathrm{euc}}
\|v\|_{A_0(\mathbb R^n)}.
\end{equation}
Then the corresponding continuous embedding holds:
\[
A_0(M,\mathbf{E})
\hookrightarrow
A_1(M,\mathbf{E}).
\]
The global constant depends on $C_{\mathrm{euc}}$, the rank of $\mathbf{E}$, and the constants $c_A,C_A$ of Proposition~\ref{prop:independencia-espacios-haz}, but not on the section.
\end{lemma}

\begin{proof}
We begin with the boundary model. The spaces on $\mathbb R^n_+$ carry the quotient norm of restriction. If $v\in A_0(\mathbb R^n_+)$ and $V\in A_0(\mathbb R^n)$ is an extension of $v$, then
\[
\|v\|_{A_1(\mathbb R^n_+)}
\leq
\|V\|_{A_1(\mathbb R^n)}
\leq
C_{\mathrm{euc}}\|V\|_{A_0(\mathbb R^n)}.
\]
Taking the infimum over all extensions proves
\begin{equation}
\label{eq:encaje-cociente-semiespacio-geometria-acotada}
\|v\|_{A_1(\mathbb R^n_+)}
\leq
C_{\mathrm{euc}}\|v\|_{A_0(\mathbb R^n_+)}.
\end{equation}
The constant is the same as on the full space.

Now fix an admissible quadratic localization system $\mathcal Q=(U_i,\phi_i,h_i)_{i\in I}$ and the synchronous frames used in Definition~\ref{def:espacios-secciones-haz-geometria-acotada}. If $\mathbf{u}=\displaystyle\sum_{a=1}^r u_i^a\mathbf{e}_{i,a}$ in $U_i$, set
\[
a_{i,a}
:=
\|(h_i u_i^a)\circ\phi_i^{-1}\|_{A_0(\mathbb X_i)},
\qquad
b_{i,a}
:=
\|(h_i u_i^a)\circ\phi_i^{-1}\|_{A_1(\mathbb X_i)}.
\]
Estimate \eqref{eq:encaje-local-modelo-geometria-acotada}, or \eqref{eq:encaje-cociente-semiespacio-geometria-acotada}, gives
$b_{i,a}\leq C_{\mathrm{euc}}a_{i,a}$
with the same constant for all indices. Since $p_0\leq p_1$, the elementary inclusion $\ell^{p_0}(I\times\{1,\dots,r\})\hookrightarrow
\ell^{p_1}(I\times\{1,\dots,r\})$ yields
\begin{align*}
\|\mathbf{u}\|_{A_1(M,\mathbf{E})}=
\left(\sum_{i\in I}\sum_{a=1}^r b_{i,a}^{p_1}\right)^{\frac{1}{p_1}}\leq
C_{\mathrm{euc}}
\left(\sum_{i\in I}\sum_{a=1}^r a_{i,a}^{p_1}\right)^{\frac{1}{p_1}}\leq
C_{\mathrm{euc}}
\left(\sum_{i\in I}\sum_{a=1}^r a_{i,a}^{p_0}\right)^{\frac{1}{p_0}}
=
C_{\mathrm{euc}}\|\mathbf{u}\|_{A_0(M,\mathbf{E})}.
\end{align*}

This calculation has been carried out in a fixed system. Its intrinsic nature follows from Proposition~\ref{prop:independencia-espacios-haz}. Recall where the geometric hypotheses enter there. On an overlap, insert the quadratic partition, multiply each component by an auxiliary cutoff equal to one on the support being transported, and compose with the coordinate change. Uniform derivative bounds for the cutoffs, transition matrices, and coordinate changes control multipliers and compositions; Jacobian bounds compare the measures. Finally, uniform multiplicity bounds both the number of summands entering a chart and the number of times a localization appears in the global sum. Thus $c_A$ and $C_A$ are also independent of $i$.
\end{proof}

\begin{theorem}[Embeddings between the real-order scales]
\label{teo:encajes-escalas-geometria-acotada}
Let $M$ be an $n$-dimensional Riemannian manifold of bounded geometry, with or without boundary, and let $\mathbf{E}\to M$ be a bundle of bounded geometry. Let $s_0,s_1\in\mathbb R$ and $1<p_0\leq p_1<\infty$.
\begin{enumerate}[label=(\alph*)]
\item If \begin{equation}
\label{eq:indices-estrictos-encajes-geometria-acotada}
s_0-\frac{n}{p_0}
>
s_1-\frac{n}{p_1},
\end{equation} then, for all $1\leq q_0,q_1\leq\infty$,
\[
F^{s_0}_{p_0,q_0}(M,\mathbf{E})
\hookrightarrow
F^{s_1}_{p_1,q_1}(M,\mathbf{E})
\]
and
\[
B^{s_0}_{p_0,q_0}(M,\mathbf{E})
\hookrightarrow
B^{s_1}_{p_1,q_1}(M,\mathbf{E}).
\]
\item If \begin{equation}
\label{eq:indices-iguales-encajes-geometria-acotada}
s_0-\frac{n}{p_0}
=
s_1-\frac{n}{p_1}
\end{equation} and $1\leq q_0\leq q_1\leq\infty$, then the same embeddings hold in the scales $F$ and $B$.
\item If \[
s_0-\frac{n}{p_0}
\geq
s_1-\frac{n}{p_1},
\] then
\[
H^{s_0,p_0}(M,\mathbf{E})
\hookrightarrow
H^{s_1,p_1}(M,\mathbf{E}).
\]
\end{enumerate}
All embeddings hold both without boundary and with boundary.
\end{theorem}

\begin{proof}
First suppose that \eqref{eq:indices-iguales-encajes-geometria-acotada}. Part (a) of Theorem~\ref{teo:encajes-escalas-orden-real}, applied with fine index $q_0$, gives
\[
F^{s_0}_{p_0,q_0}(\mathbb R^n)
\hookrightarrow
F^{s_1}_{p_1,q_0}(\mathbb R^n).
\]
Part (c) of the same theorem and $q_0\leq q_1$ then give
\[
F^{s_1}_{p_1,q_0}(\mathbb R^n)
\hookrightarrow
F^{s_1}_{p_1,q_1}(\mathbb R^n).
\]
Lemma~\ref{lem:transferencia-encaje-local-geometria-acotada} transfers the composition to $M$.

If \eqref{eq:indices-estrictos-encajes-geometria-acotada} holds, define
\[
\sigma_1
:=
s_0-\frac{n}{p_0}+\frac{n}{p_1}.
\]
Then $\sigma_1>s_1$. The result already proved, with the same fine index $q_0$, and part (b) of Theorem~\ref{teo:encajes-escalas-orden-real} yield
\[
F^{s_0}_{p_0,q_0}(M,\mathbf{E})
\hookrightarrow
F^{\sigma_1}_{p_1,q_0}(M,\mathbf{E})
\hookrightarrow
F^{s_1}_{p_1,q_1}(M,\mathbf{E}).
\]
This proves the assertions for the Triebel--Lizorkin scale. Taking $q_0=q_1=2$ and using the Euclidean identification $H^{s,p}(\mathbb R^n)=F^s_{p,2}(\mathbb R^n)$ in Theorem~\ref{teo:identificaciones-H-F-W-entero}, summed over the charts as in the proof of the lemma, proves part (c).

We turn to Besov spaces. When the indices are equal, choose $\delta>0$. The $F$ result already proved, applied at both endpoints, gives the continuous operators
\[
F^{s_0-\delta}_{p_0,p_0}(M,\mathbf{E})
\hookrightarrow
F^{s_1-\delta}_{p_1,p_1}(M,\mathbf{E}),
\qquad
F^{s_0+\delta}_{p_0,p_0}(M,\mathbf{E})
\hookrightarrow
F^{s_1+\delta}_{p_1,p_1}(M,\mathbf{E}),
\]
since $p_0\leq p_1$ is also the inequality between the fine indices of these pairs. Interpolate the identity by Theorem~\ref{teo:interpolacion-operadores-metodo-K}. The definition of the Besov scale gives
\[
B^{s_0}_{p_0,q_0}(M,\mathbf{E})
\hookrightarrow
B^{s_1}_{p_1,q_0}(M,\mathbf{E}).
\]
Monotonicity of the second index for the real method, Corollary~\ref{cor:monotonia-q-interpolacion-real}, allows $q_0$ to be replaced by $q_1\geq q_0$ in the target space. If the differential indices satisfy the strict inequality, factor once again through $B^{\sigma_1}_{p_1,q_0}(M,\mathbf{E})$. Fix $a<s_1<\sigma_1<b$ and write, using Theorem~\ref{teo:reiteracion-interpolacion-real} for independence of the auxiliary endpoints in the definition of the scale,
\[
B^{\sigma_1}_{p_1,q_0}(M,\mathbf{E})
=
(X_0,X_1)_{\theta_\sigma,q_0},
\qquad
B^{s_1}_{p_1,q_1}(M,\mathbf{E})
=
(X_0,X_1)_{\theta_s,q_1},
\]
where $X_0=F^a_{p_1,p_1}(M,\mathbf{E})$, $X_1=F^b_{p_1,p_1}(M,\mathbf{E})$, and $0<\theta_s<\theta_\sigma<1$. The $F$ result already proved gives $X_1\hookrightarrow X_0$. Lemma~\ref{lem:discretizacion-funcional-K}, separating positive and negative dyadic indices, shows that \[
(X_0,X_1)_{\theta_\sigma,q_0}
\hookrightarrow
(X_0,X_1)_{\theta_s,q_1}
\] for all $q_0,q_1$: the geometric factor $t^{\theta_\sigma-\theta_s}$ appears in the tail corresponding to $t\leq1$, whereas for $t\geq1$ we use the bound of the $X_0$ norm by the interpolation norm, proved together with completeness of the $K$ method. This completes the proof.
\end{proof}

\begin{corollary}[Integer-order Sobolev embeddings]
\label{cor:encajes-sobolev-enteros-geometria-acotada}
Let $M$ be an $n$-dimensional Riemannian manifold of bounded geometry, with or without boundary, and let $\mathbf{E}\to M$ be a bundle of bounded geometry. Let $m_0,m_1\in\mathbb N_0$ and $1<p_0\leq p_1<\infty$. If \[
m_0-\frac{n}{p_0}
\geq
m_1-\frac{n}{p_1},
\] then
\begin{equation}
\label{eq:encaje-entero-indices-geometria-acotada}
W^{m_0,p_0}(M,\mathbf{E})
\hookrightarrow
W^{m_1,p_1}(M,\mathbf{E}).
\end{equation}
In particular, let $j\in\mathbb N_0$, $m\in\mathbb N$, and $1<p<\infty$.
\begin{enumerate}[label=(\alph*)]
\item If $mp<n$ and \(
p_m^*:=\frac{np}{n-mp}
\), then
\[
W^{j+m,p}(M,\mathbf{E})
\hookrightarrow
W^{j,q}(M,\mathbf{E})
\qquad
\text{for }p\leq q\leq p_m^*.
\]
\item If $mp=n$, the same embedding holds for every $p\leq q<\infty$.
\item If $mp>n$, the same embedding holds for every $p\leq q<\infty$.
\end{enumerate}
\end{corollary}

\begin{proof}
Proposition~\ref{prop:independencia-espacios-haz} identifies $H^{m_\nu,p_\nu}(M,\mathbf{E})$ with $W^{m_\nu,p_\nu}(M,\mathbf{E})$ and provides constants $c_{m_\nu,p_\nu},C_{m_\nu,p_\nu}>0$ such that, for every $\mathbf{u}\in W^{m_\nu,p_\nu}(M,\mathbf{E})$ and $\nu\in\{0,1\}$,
\[
c_{m_\nu,p_\nu}\|\mathbf{u}\|_{W^{m_\nu,p_\nu}(M,\mathbf{E})}
\leq\|\mathbf{u}\|_{H^{m_\nu,p_\nu}(M,\mathbf{E})}
\leq C_{m_\nu,p_\nu}\|\mathbf{u}\|_{W^{m_\nu,p_\nu}(M,\mathbf{E})}.
\]
Thus \eqref{eq:encaje-entero-indices-geometria-acotada} is part (c) of Theorem~\ref{teo:encajes-escalas-geometria-acotada}. The remaining three assertions follow by writing
\[
j+m-\frac np
\geq
j-\frac nq.
\]
In the subcritical case this inequality is equivalent to $q\leq p_m^*$; in the critical and supercritical cases it holds for every finite exponent $q\geq p$.
\end{proof}

\begin{corollary}[Subcritical and critical cases]
\label{cor:encajes-Lq-geometria-acotada}
Let $M$ be an $n$-dimensional Riemannian manifold of bounded geometry, with or without boundary, and let $\mathbf{E}\to M$ be a bundle of bounded geometry. Let $s>0$ and $1<p<\infty$.
\begin{enumerate}[label=(\alph*)]
\item If $sp<n$ and \(
p_s^*:=\frac{np}{n-sp}
\), then
\[
H^{s,p}(M,\mathbf{E})
\hookrightarrow
L^q(M,\mathbf{E})
\qquad
\text{for }p\leq q\leq p_s^*.
\]
\item If $sp=n$, then
\[
H^{s,p}(M,\mathbf{E})
\hookrightarrow
L^q(M,\mathbf{E})
\qquad
\text{for every }p\leq q<\infty.
\]
\end{enumerate}
\end{corollary}

\begin{proof}
Apply part (c) of Theorem~\ref{teo:encajes-escalas-geometria-acotada} with $s_1=0$ and $p_1=q$. Proposition~\ref{prop:independencia-espacios-haz}, at order zero, identifies $H^{0,q}(M,\mathbf{E})$ with $L^q(M,\mathbf{E})$. If $sp<n$, the condition on the indices is equivalent to $q\leq p_s^*$; if $sp=n$, it is strict for every $q<\infty$.
\end{proof}

To treat the supercritical case, we need a uniform formulation of Hölder sections. In the absence of a boundary, the intrinsic norm of Definition~\ref{def:secciones-holder-haz} is already available. Near the boundary, we use the localized formulation of Definition~\ref{def:secciones-holder-haz-frontera}, replacing the finite sum by a uniform supremum.

\begin{definition}[Uniform Hölder sections up to the boundary]
\label{def:secciones-holder-uniformes-geometria-acotada}
Let $M$ be a Riemannian manifold of bounded geometry, with or without boundary, and let $\mathbf{E}\to M$ be a bundle of bounded geometry. Let $k\in\mathbb N_0$ and $0<\alpha\leq1$. Fix an admissible quadratic localization system $\mathcal Q=(U_i,\phi_i,h_i)_{i\in I}$ and synchronous frames $\mathbf e_i$. Let $\boldsymbol{\Omega}_i:=\phi_i(U_i)$ denote the bounded coordinate domain; in a boundary chart it is contained in $\mathbb R^n_+$. If \[
h_i \mathbf{u}
=
\sum_{a=1}^r u_i^a \mathbf{e}_{i,a},
\], set
\begin{equation}
\label{eq:norma-holder-uniforme-geometria-acotada}
\|\mathbf{u}\|_{\Gamma^{k,\alpha}_{\mathcal Q}(\mathbf{E})}
:=
\sup_{i\in I}
\sum_{a=1}^r
\|u_i^a\circ\phi_i^{-1}\|_{C^{k,\alpha}(\overline{\boldsymbol{\Omega}_i})}.
\end{equation}
Here we use the extension norm of $C^{k,\alpha}(\overline{\boldsymbol{\Omega}_i})$ from Definition~\ref{def:espacios-holder-euclidianos}. Define $\Gamma^{k,\alpha}(\mathbf{E})$ to be the space of sections of class $C^k$ up to the boundary for which \eqref{eq:norma-holder-uniforme-geometria-acotada} is finite.
\end{definition}

\begin{proposition}[Independence and characterization of the uniform Hölder norm]
\label{prop:holder-uniforme-geometria-acotada}
Let $M$ be a Riemannian manifold of bounded geometry, with or without boundary, and let $\mathbf{E}\to M$ be a bundle of bounded geometry. If $\mathcal Q$ and $\widetilde{\mathcal Q}$ are two admissible quadratic localization systems equipped with synchronous frames, there exist $c_{\mathcal Q,\widetilde{\mathcal Q}},
C_{\mathcal Q,\widetilde{\mathcal Q}}>0$, independent of $\mathbf{u}$, such that
\[
 c_{\mathcal Q,\widetilde{\mathcal Q}}
 \|\mathbf{u}\|_{\Gamma^{k,\alpha}_{\mathcal Q}(\mathbf{E})}
 \leq
 \|\mathbf{u}\|_{\Gamma^{k,\alpha}_{\widetilde{\mathcal Q}}(\mathbf{E})}
 \leq
 C_{\mathcal Q,\widetilde{\mathcal Q}}
 \|\mathbf{u}\|_{\Gamma^{k,\alpha}_{\mathcal Q}(\mathbf{E})}.
\]
The constants depend on $k,\alpha$, the multiplicity, and the uniform bounds on coordinate changes, cutoffs, and frames, but not on the chart indices. If $M$ has no boundary, there exist $c_{\mathrm{int}},C_{\mathrm{int}}>0$ such that
\[
 c_{\mathrm{int}}\|\mathbf{u}\|_{\Gamma^{k,\alpha}(\mathbf{E})}
 \leq\|\mathbf{u}\|_{\Gamma^{k,\alpha}_{\mathcal Q}(\mathbf{E})}
 \leq C_{\mathrm{int}}\|\mathbf{u}\|_{\Gamma^{k,\alpha}(\mathbf{E})},
\]
where the first norm is the intrinsic norm of \eqref{eq:norma-gamma-k-alpha}. If $M$ is compact with boundary and $\mathcal A$ are the finite data of Definition~\ref{def:secciones-holder-haz-frontera}, there exist $c_\partial,C_\partial>0$ such that
\[
 c_\partial\|\mathbf{u}\|_{\Gamma^{k,\alpha}_{\mathcal A}(\mathbf{E})}
 \leq\|\mathbf{u}\|_{\Gamma^{k,\alpha}_{\mathcal Q}(\mathbf{E})}
 \leq C_\partial\|\mathbf{u}\|_{\Gamma^{k,\alpha}_{\mathcal A}(\mathbf{E})}.
\]
\end{proposition}

\begin{proof}
Let $\mathcal Q=(U_i,\phi_i,h_i)_{i\in I}$ and $\widetilde{\mathcal Q}=(\widetilde U_j,\widetilde\phi_j,
\widetilde h_j)_{j\in J}$ be two admissible quadratic localization systems equipped with synchronous frames $\mathbf e_i$ and $\widetilde{\mathbf e}_j$. Insert $\displaystyle\sum_{i\in I} h_i^2=1$ into each localization of the second system:
\begin{equation}
\label{eq:transferencia-holder-uniforme-cuadratica}
\widetilde h_j \mathbf{u}
=
\sum_{\substack{i\in I\\U_i\cap\widetilde U_j\neq\varnothing}}
\widetilde h_jh_i(h_i \mathbf{u}).
\end{equation}
There are at most $L$ summands, with $L$ independent of $j$. On an overlap, the component of each summand in the frame $\widetilde{\mathbf{e}}_j$ is obtained from a component of $h_i \mathbf{u}$ by three operations: composition with the coordinate change, multiplication by an entry of the frame transition matrix, and multiplication by the two cutoffs. Also insert auxiliary cutoff functions equal to one on the supports that occur, so that all expressions are compactly supported in the exact domain of the coordinate change.

Proposition~\ref{prop:atlas-frontera-uniforme} and Proposition~\ref{prop:marcos-sincronos-uniformes} bound the derivatives of the coordinate changes, their inverses, the cutoffs, and the transition matrices uniformly in $i,j$. The Faà di Bruno formula in Theorem~\ref{faa di bruno multivariable} and the Leibniz rule then give, for each local component $f$, if $D_{ij}$ is the exact coordinate domain of the term under consideration and $D_i:=\Phi_{ij}(D_{ij})$,
\begin{equation}
\label{eq:cota-operador-holder-uniforme}
\|b_{ij}(f\circ\Phi_{ij})\|_{C^{k,\alpha}(\overline{D_{ij}})}
\leq
C_{k,\alpha}\|f\|_{C^{k,\alpha}(\overline{D_i})},
\end{equation}
where $b_{ij}$ contains the preceding multipliers and $\Phi_{ij}$ is the coordinate change. For the seminorm of order $k$, use
\[
[fg]_{C^{0,\alpha}(\overline{D_{ij}})}
\leq
\|f\|_{C^0(\overline{D_{ij}})}[g]_{C^{0,\alpha}(\overline{D_{ij}})}
+
\|g\|_{C^0(\overline{D_{ij}})}[f]_{C^{0,\alpha}(\overline{D_{ij}})}
\]
and the uniform Lipschitz bound on $\Phi_{ij}$ over the cutoff's support. In boundary charts, the coordinate changes and multipliers extend to Euclidean neighborhoods of the corresponding compact sets; thus \eqref{eq:cota-operador-holder-uniforme} descends to the extension norms upon taking the infimum over extensions.

Apply \eqref{eq:cota-operador-holder-uniforme} to \eqref{eq:transferencia-holder-uniforme-cuadratica}, sum over the $r$ components, and use the fact that there are at most $L$ terms. This yields
\[
\|\mathbf{u}\|_{\Gamma^{k,\alpha}_{\widetilde{\mathcal Q}}(\mathbf{E})}
\leq
LC_{k,\alpha}
\|\mathbf{u}\|_{\Gamma^{k,\alpha}_{\mathcal Q}(\mathbf{E})}.
\]
Interchanging the two systems proves independence.

Now suppose that $M$ has no boundary. The proof of Proposition~\ref{prop:equivalencia-norma-holder-localizada-haces}, applied in geodesic charts and synchronous frames, has uniform constants: they depend only on the constants $c_g,C_g$ in the uniform metric inequalities, the bounds on the connection coefficients, and the derivatives of the cutoffs, all independent of the chart. This gives
\[
\|\mathbf{u}\|_{\Gamma^{k,\alpha}_{\mathcal Q}(\mathbf{E})}
\leq
C\|\mathbf{u}\|_{\Gamma^{k,\alpha}(\mathbf{E})}.
\]
For the reverse inequality, use $\mathbf{u}=\displaystyle\sum_{i\in I} h_i(h_i \mathbf{u})$. Each section $h_i(h_i \mathbf{u})$ extends by zero outside $U_i$, and its intrinsic Hölder norm is bounded by the right-hand side of \eqref{eq:norma-holder-uniforme-geometria-acotada} with a uniform constant. At any point, at most $L$ summands contribute. When comparing two points, only indices whose terms are nonzero at one of them contribute, giving at most $2L$ summands. The triangle inequality for parallel transport then controls the Hölder seminorm, while the Leibniz rule controls the uniform norms of the derivatives. Thus
\[
\|\mathbf{u}\|_{\Gamma^{k,\alpha}(\mathbf{E})}
\leq
C'\|\mathbf{u}\|_{\Gamma^{k,\alpha}_{\mathcal Q}(\mathbf{E})}.
\]
If $M$ is compact with boundary and the index set has $N$ elements, the maximum norm and the sum of the local norms satisfy
\[
 \max_{i\in\{1,\ldots,N\}}a_i
 \leq \sum_{i=1}^Na_i
 \leq N\max_{i\in\{1,\ldots,N\}}a_i,
 \qquad a_i\geq0.
\]
The independence already proved allows the same data to be chosen in both definitions.
\end{proof}

\begin{theorem}[Sobolev--Morrey embeddings in bounded geometry]
\label{teo:sobolev-morrey-geometria-acotada}
Let $M$ be an $n$-dimensional Riemannian manifold of bounded geometry, with or without boundary, and let $\mathbf{E}\to M$ be a bundle of bounded geometry. Let $1<p<\infty$, $k\in\mathbb N_0$, and $0<\alpha<1$.
\begin{enumerate}[label=(\alph*)]
\item If $s>0$ and \[
k+\alpha
<
s-\frac np,
\], then
\[
H^{s,p}(M,\mathbf{E})
\hookrightarrow
\Gamma^{k,\alpha}(\mathbf{E}).
\]
\item If $m\in\mathbb N$ and \begin{equation}
\label{eq:condicion-morrey-entero-geometria-acotada}
k+\alpha
\leq
m-\frac np,
\end{equation}, then
\[
W^{m,p}(M,\mathbf{E})
\hookrightarrow
\Gamma^{k,\alpha}(\mathbf{E}).
\]
\item If $0<s<1$, $sp>n$, and $\alpha=s-\frac{n}{p}$, then
\[
B^s_{p,p}(M,\mathbf{E})
\hookrightarrow
\Gamma^{0,\alpha}(\mathbf{E}).
\]
\end{enumerate}
In each case, elements of the source space have a unique representative in the target space. If $M$ has boundary, the restriction of the representative in part (b) agrees with the trace of the section.
\end{theorem}

\begin{proof}
Fix the quadratic localization system and synchronous frames of the preceding proposition. For a section $\mathbf{u}$ of the source space, write
\[
h_i \mathbf{u}=\sum_{a=1}^r u_i^a \mathbf{e}_{i,a},
\qquad
f_i^a:=u_i^a\circ\phi_i^{-1}.
\]
Thus $f_i^a$ is the $a$th component of the already localized section $h_i \mathbf{u}$; no second factor $h_i$ is introduced. In an interior chart, Theorem~\ref{teo:encajes-sobolev-fraccionarios} gives, under the strict hypothesis in part (a),
\begin{equation}
\label{eq:morrey-local-H-geometria-acotada}
\|f_i^a\|_{C^{k,\alpha}(\overline{\boldsymbol{\Omega}_i})}
\leq
C\|f_i^a\|_{H^{s,p}(\mathbb R^n)}.
\end{equation}
The constant depends only on $n,s,p,k,\alpha$. In a boundary chart, choose an extension $F_i^a\in H^{s,p}(\mathbb R^n)$ of $f_i^a\in H^{s,p}(\mathbb R^n_+)$. The same theorem controls the Hölder norm of $F_i^a$; restricting to $\overline{\boldsymbol{\Omega}_i}$ and taking the infimum over extensions gives \eqref{eq:morrey-local-H-geometria-acotada} with the quotient norm on the right-hand side. By Proposition~\ref{prop:holder-uniforme-geometria-acotada},
\begin{align*}
\|\mathbf{u}\|_{\Gamma^{k,\alpha}(\mathbf{E})}
&\leq
C\sup_{\substack{i\in I\\a\in\{1,\ldots,r_{\mathbf E}\}}}
\|f_i^a\|_{C^{k,\alpha}(\overline{\boldsymbol{\Omega}_i})}\\
&\leq
C'
\left(
\sum_{i\in I}\sum_{a=1}^r
\|f_i^a\|_{H^{s,p}(\mathbb X_i)}^p
\right)^{\frac{1}{p}}
=
C'\|\mathbf{u}\|_{H^{s,p}(M,\mathbf{E})}.
\end{align*}
This proves (a).

For (b), use the limiting exponent in Theorem~\ref{teo:sobolev-morrey-euclidiano}. In an interior chart, $f_i^a$ is supported in a compact subset of a fixed model, and its extension by zero belongs to $W^{m,p}(\mathbb R^n)$. In a boundary chart, two different extensions are required. First extend $f_i^a$ by zero from the coordinate region to $\mathbb R^n_+$; this crosses only the lateral part, since the support of $h_i$ is separated from it, and Lemma~\ref{lem:extension-cero-media-bola-semiespacio} preserves the Sobolev norm. Then apply the operator $\mathcal E_m\colon W^{m,p}(\mathbb R^n_+)\longrightarrow W^{m,p}(\mathbb R^n)$ of Corollary~\ref{teo:extension-sobolev-semiespacio}. Finally, multiply by a fixed cutoff equal to one on the compact set where the original function may be nonzero. The result belongs to $W^{m,p}_0(B_{\mathrm{euc}}(0,R))$, agrees with $f_i^a$ on the coordinate region, and satisfies
\begin{equation}
\label{eq:extension-uniforme-morrey-frontera-geometria-acotada}
\|T_i f_i^a\|_{W^{m,p}(B_{\mathrm{euc}}(0,R))}
\leq
C\|f_i^a\|_{W^{m,p}(\mathbb R^n_+)}.
\end{equation}
The radius $R$, cutoff, and constant can be chosen independently of $i$: Proposition~\ref{prop:atlas-frontera-uniforme} uses domains of a common radius, the cutoff supports lie in a common smaller coordinate region, and the operator $\mathcal E_m$ is always the same. The Euclidean theorem and \eqref{eq:extension-uniforme-morrey-frontera-geometria-acotada} give
\[
\|f_i^a\|_{C^{k,\alpha}(\overline{\boldsymbol{\Omega}_i})}
\leq
C'\|f_i^a\|_{W^{m,p}(\mathbb X_i)}
\]
in both types of charts. Take the supremum, bound it by the $\ell^p(I)$ sum, and use the inequalities with constants $c_{m,p},C_{m,p}$ in Proposition~\ref{prop:independencia-espacios-haz}. This yields
\[
\|\mathbf{u}\|_{\Gamma^{k,\alpha}(\mathbf{E})}
\leq
C\|\mathbf{u}\|_{W^{m,p}(M,\mathbf{E})}.
\]

For (c), the Euclidean identification $B^s_{p,p}(\mathbb R^n)=W^{s,p}(\mathbb R^n)$ in Corollary~\ref{teo:identificaciones-H-W-B-F} and Lemma~\ref{teo:morrey-fraccionario-slobodeckij} give, in each chart,
\[
\|f_i^a\|_{C^{0,\alpha}(\overline{\boldsymbol{\Omega}_i})}
\leq
C\|f_i^a\|_{B^s_{p,p}(\mathbb X_i)}.
\]
On the half-space, proceed with quotient norms as in (a). The definition of $B^s_{p,p}(M,\mathbf{E})$ by real interpolation, the retraction and coretraction argument of Proposition~\ref{prop:independencia-espacios-haz}, and the Euclidean identification $B^s_{p,p}(\mathbb R^n)=F^s_{p,p}(\mathbb R^n)$ imply
\[
B^s_{p,p}(M,\mathbf{E})=F^s_{p,p}(M,\mathbf{E})
\]
and provide constants $c_{BF},C_{BF}>0$ such that
\[
c_{BF}\|\mathbf{u}\|_{F^s_{p,p}(M,\mathbf{E})}
\leq\|\mathbf{u}\|_{B^s_{p,p}(M,\mathbf{E})}
\leq C_{BF}\|\mathbf{u}\|_{F^s_{p,p}(M,\mathbf{E})}.
\]
By the localized definition of the $F$ norm, if \[
N_{\mathrm{loc}}(\mathbf{u}):=
\left(\sum_{i\in I}\sum_{a=1}^r
\|f_i^a\|_{F^s_{p,p}(\mathbb X_i)}^p\right)^{\frac{1}{p}},
\], the same constants satisfy
\[
c_{BF}N_{\mathrm{loc}}(\mathbf{u})
\leq\|\mathbf{u}\|_{B^s_{p,p}(M,\mathbf{E})}
\leq C_{BF}N_{\mathrm{loc}}(\mathbf{u}).
\]
The estimate using the supremum proves (c).

For each $i$, the representatives of $f_i^1,\dots,f_i^r$ reconstruct a section $\mathbf{s}_i$ in $U_i$ that agrees almost everywhere with $h_i \mathbf{u}$. Since the support of $h_i$ lies in a smaller coordinate region, the section $h_i\mathbf{s}_i$ extends by zero outside $U_i$. The sum \[
\widetilde{\mathbf{u}}:=\sum_{i\in I}h_i\mathbf{s}_i
\] is locally finite, defines a section of class $C^k$ up to the boundary, and agrees with $\displaystyle\sum_{i\in I}h_i^2\mathbf{u}=\mathbf{u}$ as a distribution. The preceding uniform estimates and bounded multiplicity of the cover show that $\widetilde{\mathbf{u}}\in\Gamma^{k,\alpha}(\mathbf{E})$ with the bound already obtained. If two continuous representatives agree almost everywhere, the open set where they differ must be empty, since every nonempty relative ball or half-ball has positive Riemannian measure. This proves uniqueness. In the presence of a boundary, compatibility with the trace follows by approximation in $W^{m,p}(M,\mathbf{E})$ by smooth sections: the estimate just proved gives convergence in $\Gamma^{k,\alpha}(\mathbf{E})$, while Theorem~\ref{teo:traza-haz-frontera-geometria-acotada} gives convergence of the restrictions in the Besov space on the boundary. The two distributional limits agree.
\end{proof}

\begin{proposition}[Compactness after localization]
\label{prop:compacidad-local-encajes-geometria-acotada}
Let $M$ be an $n$-dimensional Riemannian manifold of bounded geometry, with or without boundary, and let $\mathbf{E}\to M$ be a bundle of bounded geometry. Let $m_0,m_1\in\mathbb N_0$ and $1\leq p_0,p_1<\infty$, and suppose that \[
m_0>m_1,
\qquad
m_0-\frac{n}{p_0}
>
m_1-\frac{n}{p_1}.
\]. For each $\chi\in C_c^\infty(M)$, the operator \[
\mathbf{u}\longmapsto\chi \mathbf{u}
\] is compact from $W^{m_0,p_0}(M,\mathbf{E})$ to $W^{m_1,p_1}(M,\mathbf{E})$.
\end{proposition}

\begin{proof}
The support of $\chi$ meets only finitely many members of a uniformly locally finite cover. If $(\mathbf u_j)_{j\in\mathbb N}$ is bounded in $W^{m_0,p_0}(M,\mathbf{E})$, the inequalities with constants $c^W_{m_0,p_0},C^W_{m_0,p_0}$ in Proposition~\ref{prop:independencia-espacios-haz} show that, in each of these charts and for each component, the localizations of $\chi\mathbf u_j$ are bounded in a Euclidean space $W^{m_0,p_0}(\mathbb R^n)$ and supported in a fixed compact set. In an interior chart, Theorem~\ref{rellich kondrashov generalizado}, applied to a ball containing this support, provides a subsequence converging in $W^{m_1,p_1}(\mathbb R^n)$. The localized functions belong to $W^{m_0,p_0}_0(\mathbb R^n)$ of the ball, since their support is separated from its boundary. Multiplying by a cutoff equal to one on the common support and extending by zero gives convergence on all of $\mathbb R^n$. In a boundary chart, extend by zero only across the artificial faces to the half-space. The support remains in a compact subset of the closed half-space, possibly meeting $t=0$; Theorem~\ref{teo:encajes-localizados-semiespacio} gives a convergent subsequence there. Both theorems apply precisely because of the two hypotheses $m_0>m_1$ and the strict inequality between the indices, and both allow all finite exponents in the statement.

In the range $1<p_0\leq p_1<\infty$, one also obtains the following factorization through Bessel spaces. The strict inequality between the differential indices allows us to define
\[
\sigma_1
:=
m_0-\frac{n}{p_0}+\frac{n}{p_1}
>
m_1.
\]
Theorem~\ref{teo:encajes-escalas-orden-real} gives the continuous embedding $H^{m_0,p_0}(\mathbb R^n)\hookrightarrow H^{\sigma_1,p_1}(\mathbb R^n)$, while Theorem~\ref{teo:rellich-orden-real}, applied on a bounded open set containing the fixed support, gives the compact embedding $H^{\sigma_1,p_1}(\mathbb R^n)\hookrightarrow H^{m_1,p_1}(\mathbb R^n)$. The integer-order identifications thus provide a subsequence converging in $W^{m_1,p_1}(\mathbb R^n)$ in an interior chart. Since all functions are supported in a fixed compact set, multiplying by a cutoff equal to one there and extending by zero turns this local convergence into convergence in the full Euclidean space. In a boundary chart, first apply lateral extension by zero and then the half-space operator, as in \eqref{eq:extension-uniforme-morrey-frontera-geometria-acotada}; the same factorization argument applies to the extensions supported in a fixed compact set.

Only finitely many charts and components are involved. Successive extractions yield a single subsequence converging in all of them. The synthesis estimate, which uses multipliers, coordinate changes, and uniform multiplicity, reconstructs convergence in $W^{m_1,p_1}(M,\mathbf{E})$. Thus the operator is compact.
\end{proof}

\begin{proposition}[Interpolation inequality with a small parameter]
\label{prop:ehrling-geometria-acotada}
Let $M$ be an $n$-dimensional Riemannian manifold of bounded geometry, with or without boundary, and let $\mathbf{E}\to M$ be a bundle of bounded geometry. Let $t,s,s_0\in\mathbb N_0$ and $p_0,p,q\in(1,\infty)$ satisfy \[
t<s<s_0,
\qquad
\max\{p_0,q\}\leq p,
\] and \begin{equation}
\label{eq:indices-ehrling-geometria-acotada}
s-\frac np
<
s_0-\frac n{p_0}.
\end{equation} For every $\varepsilon>0$ there exists $C_\varepsilon>0$ such that
\begin{equation}
\label{eq:ehrling-geometria-acotada}
\|\mathbf{u}\|_{W^{s,p}(M,\mathbf{E})}
\leq
\varepsilon\|\mathbf{u}\|_{W^{s_0,p_0}(M,\mathbf{E})}
+
C_\varepsilon\|\mathbf{u}\|_{W^{t,q}(M,\mathbf{E})}
\end{equation}
for every $\mathbf{u}\in W^{s_0,p_0}(M,\mathbf{E})\cap W^{t,q}(M,\mathbf{E})$.
\end{proposition}

\begin{proof}
First establish an inequality on the fixed local model. Let $K$ denote the common compact set containing the supports of the localization functions. The strict inequality in \eqref{eq:indices-ehrling-geometria-acotada} allows us to set
\[
\sigma
:=
s_0-\frac{n}{p_0}+\frac{n}{p}
>
s.
\]
Theorem~\ref{teo:encajes-escalas-orden-real} gives \[
H^{s_0,p_0}(\mathbb R^n)
\hookrightarrow
H^{\sigma,p}(\mathbb R^n).
\]. Choose a bounded Lipschitz domain $\Omega$ and a cutoff $\eta\in C_c^\infty(\Omega)$ equal to one on a neighborhood of $K$. Theorem~\ref{teo:rellich-orden-real} shows that the restriction of a bounded family in $H^{\sigma,p}(\mathbb R^n)$ has a subsequence converging in $H^{s,p}(\Omega)$. Since each function under consideration is supported in $K$, it agrees on the whole space with the extension by zero of $\eta v\restriction_\Omega$. Continuity of multiplication by $\eta$ and of this extension, whose support remains separated from $\partial\Omega$, proves that \[
\{v\in W^{s_0,p_0}(\mathbb R^n)\mid\operatorname{supp}v\subseteq K\}
\hookrightarrow
W^{s,p}(\mathbb R^n)
\] is compact. Here the identifications $H^{s_0,p_0}(\mathbb R^n)=W^{s_0,p_0}(\mathbb R^n)$ and $H^{s,p}(\mathbb R^n)=W^{s,p}(\mathbb R^n)$ are at integer orders. On the other hand, since $t<s$ and $q\leq p$, first lose regularity while keeping the exponent $p$ fixed and then apply Hölder's inequality on the bounded domain containing $K$; this gives the continuous embedding into $W^{t,q}(\mathbb R^n)$ of the subspace of functions in $W^{s,p}(\mathbb R^n)$ supported in $K$.

We claim that, for each $\delta>0$, there exists $C_\delta>0$ such that
\begin{equation}
\label{eq:ehrling-local-modelo-geometria-acotada}
\|v\|_{W^{s,p}(\mathbb R^n)}
\leq
\delta\|v\|_{W^{s_0,p_0}(\mathbb R^n)}
+
C_\delta\|v\|_{W^{t,q}(\mathbb R^n)}
\end{equation}
for every $v\in W^{s_0,p_0}(\mathbb R^n)$ supported in $K$. If this failed for some $\delta$, there would be a sequence $(v_j)\subset W^{s_0,p_0}(\mathbb R^n)$ supported in $K$ which, after normalization, would satisfy
\[
\|v_j\|_{W^{s,p}(\mathbb R^n)}=1,
\qquad
\|v_j\|_{W^{s_0,p_0}(\mathbb R^n)}\leq\delta^{-1},
\qquad
\|v_j\|_{W^{t,q}(\mathbb R^n)}\longrightarrow0.
\]
By compactness, a subsequence would converge in $W^{s,p}(\mathbb R^n)$ to a limit of norm one. By continuity of $W^{s,p}_K(\mathbb R^n)\hookrightarrow W^{t,q}(\mathbb R^n)$, the same limit would be zero in $W^{t,q}(\mathbb R^n)$ and hence the zero distribution. This contradiction proves \eqref{eq:ehrling-local-modelo-geometria-acotada}.

In a boundary chart, the compact inclusion used in the contradiction is obtained by first extending by zero from the half-ball to the half-space, only across the lateral part, then applying $\mathcal E_{s_0}$ and multiplying by a fixed cutoff. The continuous inclusion of the space of order $s$ into that of order $t$ is proved directly on the half-space by Hölder's inequality, since all functions are supported in $K$. Thus the same contradiction gives \eqref{eq:ehrling-local-modelo-geometria-acotada} with half-space norms. Since the domains, compact set $K$, cutoffs, and extension operator are common to the entire family, $C_\delta$ is independent of the chart. Apply the inequality to each localized component and take the $\ell^p(I\times\{1,\ldots,r_{\mathbf E}\})$ norm. The triangle inequality in $\ell^p(I\times\{1,\ldots,r_{\mathbf E}\})$ gives
\[
\|(b_{i,a})\|_{\ell^p(I\times\{1,\ldots,r_{\mathbf E}\})}
\leq
\delta\|(a_{i,a})\|_{\ell^p(I\times\{1,\ldots,r_{\mathbf E}\})}
+
C_\delta\|(c_{i,a})\|_{\ell^p(I\times\{1,\ldots,r_{\mathbf E}\})},
\]
where $a,b,c$ are the local norms of orders and exponents $(s_0,p_0)$, $(s,p)$, and $(t,q)$, respectively. The hypotheses $p_0\leq p$ and $q\leq p$ imply
\[
\|(a_{i,a})\|_{\ell^p(I\times\{1,\ldots,r_{\mathbf E}\})}
\leq
\|(a_{i,a})\|_{\ell^{p_0}(I\times\{1,\ldots,r_{\mathbf E}\})},
\qquad
\|(c_{i,a})\|_{\ell^p(I\times\{1,\ldots,r_{\mathbf E}\})}
\leq
\|(c_{i,a})\|_{\ell^q(I\times\{1,\ldots,r_{\mathbf E}\})}.
\]
Proposition~\ref{prop:independencia-espacios-haz} compares these three sums with the global covariant norms. Choose $\delta$ so that the first coefficient, after incorporating the constants $c_A,C_A$ of Proposition~\ref{prop:independencia-espacios-haz}, is $\varepsilon$. This gives \eqref{eq:ehrling-geometria-acotada}.
\end{proof}

\begin{remark}[Local compactness and loss of global compactness]
\label{obs:perdida-compacidad-global-geometria-acotada}
Proposition~\ref{prop:compacidad-local-encajes-geometria-acotada} does not imply compactness of the global embedding. Even on $M=\mathbb R^n$, if $0\neq v\in C_c^\infty(\mathbb R^n)$ and $v_j(x):=v(x-je_1)$, the sequence $(v_j)$ is bounded in all Sobolev spaces of a fixed order. Choosing a subsequence with disjoint supports gives, for $j\neq\ell$,
\[
\|v_j-v_\ell\|_{W^{m,p}(\mathbb R^n)}^p
=
2\|v\|_{W^{m,p}(\mathbb R^n)}^p.
\]
Thus there is no convergent subsequence. Bounded geometry prevents the local constants from degenerating, but does not prevent a sequence from escaping to infinity. Inequality \eqref{eq:ehrling-geometria-acotada} is compatible with this phenomenon: it is a global estimate with a small parameter, not a statement of global compactness.
\end{remark}

\begin{remark}[Completeness and weighted embeddings]
\label{obs:encajes-ponderados-variedades-completas}
Completeness alone does not provide a uniform coordinate radius or the preceding unweighted embeddings. Fix $0<\varepsilon<1$. A ball $B_{\mathbf{g}}(x,R)$ is called $(1,\varepsilon)$-admissible if it has a chart $\phi\colon B_{\mathbf{g}}(x,R)\longrightarrow\Omega\subseteq\mathbb R^n$ whose metric coefficients $g_{ij}(z)$ satisfy
\[
 (1-\varepsilon)\sum_{i=1}^n(\xi^i)^2
 \leq\sum_{i,j=1}^n g_{ij}(z)\xi^i\xi^j
 \leq(1+\varepsilon)\sum_{i=1}^n(\xi^i)^2
 \qquad(z\in\Omega,\ \xi\in\mathbb R^n)
\]
and
\[
R\max_{\substack{\beta\in\mathbb N_0^n\\|\beta|=1}}\max_{i,j\in\{1,\ldots,n\}}
\|\partial^\beta g_{ij}\|_{L^\infty(\phi(B_{\mathbf{g}}(x,R)))}
\leq\varepsilon.
\]
In the first condition, $\xi$ represents the components of a tangent vector in the chart; hence this condition is equivalent to
\[
 \sqrt{1-\varepsilon}\,\|d\phi_p(v)\|
 \leq|v|_{\mathbf g}
 \leq\sqrt{1+\varepsilon}\,\|d\phi_p(v)\|
 \qquad(p\in B_{\mathbf g}(x,R),\ v\in T_pM).
\]
The second condition additionally controls the first coordinate derivatives of the metric. If $R'(x)$ is the supremum of these radii, set $\displaystyle R_\varepsilon(x):=\min\left\{1,\frac{R'(x)}{2}\right\}$. The local constants can be retained inside the integral through weights formed from powers of $R_\varepsilon$. With the convention \[
\|\mathbf{u}\|_{L^a(M,\mathbf{E};w)}
:=
\left(\int_M |\mathbf{u}|_{\mathbf{h}_{\mathbf{E}}}^a w\,d\lambda_{\mathbf{g}}\right)^{\frac{1}{a}},
\] for $1<r<n$ and $r^*=\frac{nr}{n-r}$, one obtains an estimate of the form
\[
\|\mathbf{u}\|_{L^{r^*}(M,\mathbf{E};R_\varepsilon^{2r^*})}
\leq
C\|\mathbf{u}\|_{W^{1,r}(M,\mathbf{E})}
\]
on a complete manifold, under the local control of the connection needed for sections; see \cite[Theorem~7.5]{Amar2019SobolevWeights}. When $R_\varepsilon$ is uniformly bounded away from zero, there exists $\rho>0$ such that $\rho\leq R_\varepsilon(x)\leq1$ for every $x\in M$, and hence
\[
\rho^{2r^*}\leq R_\varepsilon(x)^{2r^*}\leq1.
\]
These inequalities recover the unweighted mechanism of this section.
\end{remark}

\part{Index theory}
\chapter*{Introduction to Part IV}
\addcontentsline{toc}{chapter}{Introduction to Part IV}
\markboth{Part IV. Index theory}{Introduction to Part IV}

One of the deepest questions in elliptic theory is how much can be known about an operator without knowing all its coefficients. For an elliptic operator on a closed manifold, lower-order terms influence the particular solutions but do not alter certain fundamental global information. Among this information is the index, defined as the difference between the dimension of the kernel and the dimension of the cokernel.

The first step will be to extend the differential calculus to the pseudodifferential calculus. Its symbols allow parametrices to be constructed and the regularity and singularities of solutions to be described precisely. Fredholm theory turns a parametrix modulo compact operators into a stable integer. This integer still has an analytic origin, but the principal symbol has an interpretation independent of the full operator: it determines a $K$-theory class on the cotangent bundle.

The second half of this part develops the topological tools needed to interpret that class. Bott periodicity, the Thom isomorphism, and the morphisms associated with embeddings allow the problem to be transported to Euclidean space. The Atiyah--Singer theorem finally asserts that the index obtained by this procedure agrees with the analytic index. The question uniting these chapters is simple to state, although its answer requires several languages: how can a differential equation retain, in its leading term, topological information about the entire manifold?

\begin{semblanzaHistorica}{From elliptic analysis to the index}
Fredholm discovered that certain integral equations on function spaces retained an essentially finite-dimensional structure. With the development of elliptic theory and the pseudodifferential calculus, a similar idea emerged for differential operators. In the 1960s, Atiyah and Singer proved that the analytic index could be calculated from topological data of the principal symbol. The result brought together ideas from Hodge, Hirzebruch, Thom, Bott, Grothendieck, and operator theory. Since then, the index has become one of the most fruitful bridges between analysis, geometry, and topology.
\end{semblanzaHistorica}

\chapter{Pseudodifferential operators}
\label{cap:operadores-pseudodiferenciales}

The starting point of the pseudodifferential calculus is a difficulty that already arises when one tries to invert a constant-coefficient elliptic operator. If \[
P(D)=\sum_{|\alpha|\leq m}c_\alpha D^\alpha,
\], the Fourier transform turns the equation $P(D)u=f$ into $p(\xi)\widehat u(\xi)=\widehat f(\xi)$. When $p(\xi)$ is invertible, the formal solution is obtained by multiplication by $p(\xi)^{-1}$; this inverse, however, is almost never a polynomial and therefore does not correspond to a differential operator. Variable coefficients introduce a further difficulty: the symbol depends on both the point $x$ and the frequency $\xi$, and composition no longer reduces to the ordinary product of functions.

Pseudodifferential operators enlarge the class of differential operators just enough to solve this problem. Their symbols allow separate control of variation in the base and behavior at high frequencies. Composition is described by an asymptotic expansion, the adjoint remains in the same class, and ellipticity leads to an inverse modulo operators with smooth kernels. This last property turns symbolic inversion into elliptic regularity and, on a closed manifold, into a Fredholm statement.

The analytic ingredients were prepared in earlier chapters. The Fourier transform and tempered distributions of Chapter~\ref{cap:distribuciones-fourier} provide the local model; the Bessel potentials of Chapter~\ref{cap:regularidad-intermedia} show that natural multipliers need not be polynomial; the principal symbol and ellipticity of Chapter~\ref{cap:simbolo-principal-operadores-elipticos} identify the part to be inverted. For example, \[
J^s=(I-\Delta)^{\frac{s}{2}},\qquad
\widehat{J^su}(\xi)=\langle\xi\rangle^s\widehat u(\xi),
\] is an operator of order $s$, whereas $J^{-2}=(I-\Delta)^{-1}$ gains two derivatives without being differential. The calculus we construct explains this phenomenon within a stable theory and then transfers it to manifolds and vector bundles.

We begin on $\mathbb R^n$ with symbol classes, oscillatory integrals, and Kohn--Nirenberg quantization. The study of kernels will distinguish locality, pseudolocality, and proper support, and yield the composition and adjoint formulas. We then globalize the calculus, identify the kernels as distributions conormal to the diagonal, and compare local quantization with a Riemannian quantization based on the exponential map. Parametrices will then lead to microlocal regularity and the wavefront set.

The last part of the chapter separates two global situations. On manifolds of bounded geometry, we construct a uniform calculus with finite propagation that preserves Sobolev estimates over an infinite cover. On a closed manifold, smoothing remainders are compact and the parametrix yields Fredholm operators. Finally, the parameter-dependent calculus will allow us to construct Seeley's complex powers, compare them with the spectral calculus, and relate the heat kernel to the McKean--Singer formula.

The symbolic and microlocal theory will rely primarily on Hörmander~\cite{HormanderIII}, Taylor~\cite{Taylor1974Pseudo}, and Shubin~\cite{Shubin2001}. For geometric quantization we follow the idea of Bokobza--Haggiag~\cite{BokobzaHaggiag1969}; the construction of complex powers is due to Seeley~\cite{Seeley1967}, and the transition to the analytic index uses the exposition of Bleecker and Booß--Bavnbek~\cite{BleeckerBoossIndex}. The proof of the Atiyah--Singer theorem is left to later chapters.

\begin{semblanzaHistorica}{Hörmander, Seeley, and the maturation of the symbolic calculus}
The pseudodifferential calculus emerged by enlarging Fourier multipliers to a class stable under composition and capable of constructing approximate inverses for elliptic operators. Hörmander organized microlocal analysis around symbols, kernels, and wavefront sets; Seeley showed how the same language constructs complex powers of elliptic operators. The strength of the theory lies in its simultaneous treatment of three scales: local behavior in the base variable, homogeneity in frequency, and global regularity of solutions.
\end{semblanzaHistorica}

\section{From the Fourier inverse to symbol classes}
\label{sec:pseudo-simbolos-euclidianos}

\subsection{Conventions and the first model}

With the unitary normalization of the Fourier transform fixed in Definition~\ref{def:transformada-fourier-L1},
\[
\widehat u(\xi)=(2\pi)^{-\frac{n}{2}}\int_{\mathbb R^n}e^{-ix\cdot\xi}u(x)\,dx,
\qquad
u(x)=(2\pi)^{-\frac{n}{2}}\int_{\mathbb R^n}e^{ix\cdot\xi}\widehat u(\xi)\,d\xi.
\]
Moreover, $D_j$ denotes the partial derivative $\partial_j$, not the operator $-i\partial_j$. By Proposition~\ref{prop:calculo-fourier-schwartz}, $\mathcal F(D^\alpha u)(\xi)=(i\xi)^\alpha\widehat u(\xi)$. These conventions determine all factors that will appear in quantization.

First consider $P(D)=\displaystyle\sum_{|\alpha|\leq m}c_\alpha D^\alpha$. If $u,f\in\mathcal S(\mathbb R^n)$ and $P(D)u=f$, then
\[
p(\xi)\widehat u(\xi)=\widehat f(\xi),
\qquad
p(\xi):=\displaystyle\sum_{|\alpha|\leq m}c_\alpha(i\xi)^\alpha.
\]
When $p(\xi)\neq0$, formal inversion uses $p(\xi)^{-1}$. The polynomial $p$ comes from a differential operator, but its inverse is not polynomial. If the coefficients now depend on $x$, the function $p(x,\xi)^{-1}$ depends simultaneously on position and frequency. The symbol class must control derivatives in both variables while allowing each derivative with respect to $\xi$ to lower the order.

We write $\langle\xi\rangle=(1+\|\xi\|^2)^{\displaystyle\frac{1}{2}}$. All real bundles and vector spaces will be complexified during the Fourier calculus; at the end, operators preserving the real structure will be restricted back to it.

\subsection{Symbols of type \texorpdfstring{$(1,0)$}{(1,0)}}

The preceding heuristic suggests the two estimates a symbol must satisfy. Derivatives in $x$ measure variation of the coefficients, while each derivative in $\xi$ should lower the order. The class of type $(1,0)$ formalizes both controls simultaneously.

\begin{definition}[Symbol class]
\label{def:clase-simbolos-10}
\index{symbol@symbol!symbol class}
Let $n,r,s\in\mathbb N$ and $m\in\mathbb R$. The global uniform class \[
S^m_{1,0}(\mathbb R^n;\operatorname{Hom}(\mathbb C^r,\mathbb C^s))
\] consists of smooth maps $a\colon\mathbb R_x^n\times\mathbb R_\xi^n\longrightarrow\operatorname{Hom}(\mathbb C^r,\mathbb C^s)$ such that, for all multi-indices $\alpha,\beta$, there exists $C_{\alpha,\beta}>0$ with
\[
\|D_x^\alpha D_\xi^\beta a(x,\xi)\|
\leq
C_{\alpha,\beta}\langle\xi\rangle^{m-|\beta|}
\]
for every $(x,\xi)\in\mathbb R^n\times\mathbb R^n$. If $r=s=1$, we simply write $S^m_{1,0}(\mathbb R^n)$.

For an open set $U\subseteq\mathbb R^n$, the local class \[
S^m_{1,0}\!\left(
U\times\mathbb R^n;
\operatorname{Hom}(\mathbb C^r,\mathbb C^s)
\right)
\] consists of smooth maps with the same target such that, for every compact set $K\Subset U$ and all multi-indices $\alpha,\beta$, there exists $C_{K,\alpha,\beta}>0$ satisfying
\[
\|D_x^\alpha D_\xi^\beta a(x,\xi)\|
\leq C_{K,\alpha,\beta}\langle\xi\rangle^{m-|\beta|},
\qquad (x,\xi)\in K\times\mathbb R^n.
\]
When $U\neq\mathbb R^n$, we may abbreviate this local class as
$S^m_{1,0}(U;\operatorname{Hom}(\mathbb C^r,\mathbb C^s))$.
For $U=\mathbb R^n$, the local class is always written with both factors,
$S^m_{1,0}(\mathbb R^n\times\mathbb R^n;\operatorname{Hom}(\mathbb C^r,\mathbb C^s))$;
the notation with the single argument $\mathbb R^n$ remains reserved for
global uniform estimates. We omit the target only in the scalar case. When the order and all spaces are fixed within a block of text, the subscript $(1,0)$ may be abbreviated.
\end{definition}

For $n,r,s\in\mathbb N$ and $m\in\mathbb R$, the natural topology on $S^m_{1,0}(\mathbb R^n;
\operatorname{Hom}(\mathbb C^r,\mathbb C^s))$ is the Fréchet topology defined by the seminorms
\[
q_{N,m}(a)
:=
\max_{|\alpha|+|\beta|\leq N}
\sup_{(x,\xi)\in\mathbb R^{2n}}
\langle\xi\rangle^{-m+|\beta|}
\|D_x^\alpha D_\xi^\beta a(x,\xi)\|,
\qquad N\in\mathbb N_0.
\]
On an open set $U$, use a compact exhaustion $(K_j)$ and the seminorms obtained by restricting $x$ to $K_j$. For $U\subseteq\mathbb R^n$ open, define the class of \textbf{smoothing symbols} by
\index{smoothing symbol@smoothing symbol}
\begin{equation}
\begin{split}
&S^{-\infty}\!\left(
U\times\mathbb R^n;
\operatorname{Hom}(\mathbb C^r,\mathbb C^s)
\right)\\
&\qquad:=
\bigcap_{\mu\in\mathbb R}
S^\mu_{1,0}\!\left(
U\times\mathbb R^n;
\operatorname{Hom}(\mathbb C^r,\mathbb C^s)
\right).
\end{split}
\label{eq:def-simbolos-regularizantes-euclidianos}
\end{equation}
Equivalently, $a$ belongs to this intersection if, for every compact set $K\Subset U$, all multi-indices $\alpha,\beta$, and every $N\in\mathbb N_0$, there exists $C_{K,\alpha,\beta,N}>0$ such that
\begin{equation}
\|D_x^\alpha D_\xi^\beta a(x,\xi)\|
\leq C_{K,\alpha,\beta,N}\langle\xi\rangle^{-N},
\qquad (x,\xi)\in K\times\mathbb R^n.
\label{eq:estimaciones-simbolo-regularizante-euclidiano}
\end{equation}
For the global uniform class, the same estimate is required for every $x\in\mathbb R^n$, with a constant independent of $x$, and we write
\begin{equation}
S^{-\infty}(\mathbb R^n;
\operatorname{Hom}(\mathbb C^r,\mathbb C^s))
:=
\bigcap_{\mu\in\mathbb R}
S^\mu_{1,0}(\mathbb R^n;
\operatorname{Hom}(\mathbb C^r,\mathbb C^s)).
\label{eq:def-simbolos-regularizantes-globales}
\end{equation}
In any of the global or local spaces just defined, the inclusions from order $m_1$ into order $m_2$ are continuous when $m_1\leq m_2$, differentiation $D_x^\alpha D_\xi^\beta$ takes order $m$ to $m-|\beta|$, and matrix multiplication adds orders whenever the target spaces allow composition.

The reason for privileging $(\rho,\delta)=(1,0)$ can already be seen in these three properties: frequency differentiation lowers the order, position differentiation does not increase it, and the Taylor expansions appearing in composition produce strictly decreasing orders. The general classes $S^m_{\rho,\delta}(U\times\mathbb R^n;
\operatorname{Hom}(\mathbb C^r,\mathbb C^s))$ replace $m-|\beta|$ by $m-\rho|\beta|+\delta|\alpha|$; they will not be needed here.

\begin{definition}[Euclidean elliptic symbol]
\label{def:simbolo-eliptico-euclidiano}
\index{elliptic symbol@elliptic symbol}
Let $U\subseteq\mathbb R^n$ be open, and let $r\in\mathbb N$ and $m\in\mathbb R$. A matrix-valued symbol $a\in S^m_{1,0}(U;\operatorname{End}(\mathbb C^r))$ is \textbf{elliptic} if, for every compact set $K\Subset U$, there exist $R_K,C_K>0$ such that $a(x,\xi)$ is invertible and
\[
\|a(x,\xi)^{-1}\|
\leq
C_K\langle\xi\rangle^{-m}
\]
when $x\in K$ and $\|\xi\|\geq R_K$. It is \textbf{uniformly elliptic} if $R_K$ and $C_K$ can be chosen independently of $K$.
\end{definition}

\subsection{Classical symbols and asymptotic expansions}

The class $S^m_{1,0}$ controls the size of a symbol but does not distinguish successive contributions from its different degrees. To construct a parametrix, we first need to separate the leading term and correct it order by order. Classical symbols provide precisely this homogeneous expansion at high frequencies.

\begin{definition}[Classical symbol]
\label{def:simbolo-clasico-pseudo}
\index{classical symbol@classical symbol}
Let $U\subseteq\mathbb R^n$ be open, and let $r,s\in\mathbb N$ and $m\in\mathbb R$. A symbol $a\in S^m_{1,0}(U;\operatorname{Hom}(\mathbb C^r,\mathbb C^s))$ is \textbf{classical} or \textbf{polyhomogeneous} if there exist smooth maps \[
a_{m-j}\colon U\times(\mathbb R^n\setminus\{0\})\longrightarrow
\operatorname{Hom}(\mathbb C^r,\mathbb C^s),
\qquad j\in\mathbb N_0,
\], homogeneous in the sense that
$a_{m-j}(x,t\xi)=t^{m-j}a_{m-j}(x,\xi)$
for $t>0$, and a function $\vartheta\in C^\infty(\mathbb R^n)$ that vanishes near $0$ and equals one for $\|\xi\|\geq1$, such that, for every $N\geq1$,
\[
a-\displaystyle\sum_{j=0}^{N-1}\vartheta a_{m-j}
\in S^{m-N}_{1,0}(U;\operatorname{Hom}(\mathbb C^r,\mathbb C^s)).
\]
We write $a\in S^m_{\mathrm{cl}}(U;\operatorname{Hom}(\mathbb C^r,\mathbb C^s))$ and call $a_m$ the \textbf{principal homogeneous component}.
\end{definition}

The presence of $\vartheta$ does not alter the homogeneous information: two cutoffs yield symbols whose difference is compactly supported in $\xi$ and belongs to $S^{-\infty}(U\times\mathbb R^n;
\operatorname{Hom}(\mathbb C^r,\mathbb C^s))$. The component $a_m$ is unique. If $b_m$ is another principal component, then $a_m-b_m\in S^{m-1}_{1,0}(U;
\operatorname{Hom}(\mathbb C^r,\mathbb C^s))$ for $\|\xi\|\geq1$. For $\omega\in\mathbb S^{n-1}$ and $t\geq1$, homogeneity gives
\[
t^m\|a_m(x,\omega)-b_m(x,\omega)\|
\leq C_Kt^{m-1}.
\]
Dividing by $t^m$ and letting $t\to\infty$ gives $a_m=b_m$. The same inductive argument proves uniqueness of all homogeneous components.

\begin{definition}[Asymptotic expansion]
\label{def:desarrollo-asintotico-simbolos}
\index{asymptotic expansion@asymptotic expansion!of symbols}
Let $U\subseteq\mathbb R^n$ be open, let $r,s\in\mathbb N$, and let $m_0\geq m_1\geq\cdots$ be real numbers with $m_j\to-\infty$ and
\[
a_j\in S^{m_j}_{1,0}
(U;\operatorname{Hom}(\mathbb C^r,\mathbb C^s)).
\]
For $a\in S^{m_0}_{1,0}(U;
\operatorname{Hom}(\mathbb C^r,\mathbb C^s))$, we write \[
a\sim\displaystyle\sum_{j=0}^\infty a_j
\] if, for every $N\geq1$,
\[
a-\displaystyle\sum_{j=0}^{N-1}a_j
\in S^{m_N}_{1,0}(U;
\operatorname{Hom}(\mathbb C^r,\mathbb C^s)).
\]
When several $m_j$ coincide, the conclusion is understood after grouping terms of equal order.
\end{definition}

Differentiation may be performed term by term: if $a\sim\displaystyle\sum_{j=0}^{\infty} a_j$, then $D_x^\alpha D_\xi^\beta a\sim\displaystyle\sum_{j=0}^{\infty}D_x^\alpha D_\xi^\beta a_j$, with orders $m_j-|\beta|$. The coefficients of a general expansion are not unique as functions; what is determined is the class of the term of order $m_j$ modulo symbols of the next order. For classical symbols, homogeneity removes this ambiguity.

\begin{lemma}[Asymptotic summation of symbols]
\label{lem:suma-asintotica-simbolos}
\index{symbolic Borel lemma@symbolic Borel lemma}
Let $U\subseteq\mathbb R^n$ be open and let $r,s\in\mathbb N$. In the setting of Definition~\ref{def:desarrollo-asintotico-simbolos}, there exists \[
a\in S^{m_0}_{1,0}(U;
\operatorname{Hom}(\mathbb C^r,\mathbb C^s))
\] with $a\sim\displaystyle\sum_{j=0}^\infty a_j$. In particular, if $h_{m_0-j}$ is smooth on $U\times(\mathbb R^n\setminus\{0\})$ with values in $\operatorname{Hom}(\mathbb C^r,\mathbb C^s)$ and homogeneous of degree $m_0-j$ for each $j\in\mathbb N_0$, then, for a radial cutoff $\vartheta$ as in Definition~\ref{def:simbolo-clasico-pseudo}, one can choose $a\in S^{m_0}_{\mathrm{cl}}(U;
\operatorname{Hom}(\mathbb C^r,\mathbb C^s))$ with
\[
a\sim\displaystyle\sum_{j=0}^{\infty}\vartheta h_{m_0-j}.
\]
\end{lemma}

The lattice of degrees $m_0-j$ in the last clause is essential to the notion in Definition~\ref{def:simbolo-clasico-pseudo}. A homogeneous family of merely decreasing real degrees admits a general asymptotic sum, but defines only a homogeneous expansion with that set of degrees; it is classical in the sense adopted here when those degrees belong to $m_0-\mathbb N_0$, inserting zero components at the missing degrees.

\begin{proof}
Throughout this proof, $U,r,s$ remain fixed, and we write
\[
S^q:=S^q_{1,0}
\bigl(U\times\mathbb R^n;
\operatorname{Hom}(\mathbb C^r,\mathbb C^s)\bigr).
\]
Since $m_j\to-\infty$ and the sequence is nonincreasing, each order can occur only finitely many times. For each value $\mu$ that occurs, the set \[
I_\mu:=\{j\in\mathbb N_0\mid m_j=\mu\}
\] is finite. Group the terms corresponding to $I_\mu$ and replace the group by the finite sum $\displaystyle \sum_{j\in I_\mu}a_j$. This is exactly the convention fixed after Definition~\ref{def:desarrollo-asintotico-simbolos}. It suffices to treat the case
\[
m_0>m_1>m_2>\cdots.
\]

Fix an exhaustion $K_1\Subset K_2\Subset\cdots$ of $U$ and a function $\chi\in C^\infty(\mathbb R^n)$ vanishing on $\{\|\xi\|\leq1\}$ and equal to one on $\{\|\xi\|\geq2\}$. For a compact set $K\Subset U$, write
\[
q^{K}_{L,\mu}(c)
:=
\max_{|\alpha|+|\beta|\leq L}
\sup_{x\in K,\,\xi\in\mathbb R^n}
\langle\xi\rangle^{-\mu+|\beta|}
\|D_x^\alpha D_\xi^\beta c(x,\xi)\|.
\]
We will choose a sequence $\varepsilon_j\to 0^{+}$ and set
\[
b_j(x,\xi):=\chi(\varepsilon_j\xi)a_j(x,\xi).
\]
If $N<j$, then $m_j<m_N$. For $0<\varepsilon\leq1$, the
Leibniz rule gives
\[
 D_x^\alpha D_\xi^\beta(\chi(\varepsilon\xi)a_j)
 =\sum_{\gamma\leq\beta}\binom{\beta}{\gamma}
   \varepsilon^{|\gamma|}(D^\gamma\chi)(\varepsilon\xi)
   D_x^\alpha D_\xi^{\beta-\gamma}a_j.
\]
Each term is supported where $\|\xi\|\geq\varepsilon^{-1}$.
When $|\gamma|>0$, we also have
$1\leq\|\varepsilon\xi\|\leq2$; hence
\[
 \varepsilon^{|\gamma|}\langle\xi\rangle^{|\gamma|}
 \leq5^{\frac{|\gamma|}{2}},
 \qquad
 \langle\xi\rangle^{m_j-m_N}
 \leq\varepsilon^{m_N-m_j}.
\]
The first inequality also holds for $\gamma=0$. Multiplying the
Leibniz formula by the weight defining $q^K_{L,m_N}$ and using the
seminorm $q^K_{L,m_j}(a_j)$ gives
\begin{equation}
\label{eq:pequenez-corte-suma-asintotica}
 q^K_{L,m_N}\bigl(\chi(\varepsilon\,\cdot)a_j\bigr)
 \leq C_{L,\chi}\varepsilon^{m_N-m_j}q^K_{L,m_j}(a_j)
 \longrightarrow0
 \qquad(\varepsilon\longrightarrow0^+).
\end{equation}
The constant $C_{L,\chi}$ depends only on $L$ and derivatives of the
cutoff up to that order. The limit follows from $m_N-m_j>0$ and
finiteness of $q^K_{L,m_j}(a_j)$.

For each $j\geq1$, there are only finitely many conditions
\[
 q^{K_p}_{p,m_N}(b_j)\leq2^{-j},
 \qquad0\leq N<j,\qquad1\leq p\leq j.
\]
Choose $0<\varepsilon_0\leq1$. Convergence
\eqref{eq:pequenez-corte-suma-asintotica} allows successive choices
of $\varepsilon_j$ satisfying all the conditions and also
$0<\varepsilon_j<\min\{\varepsilon_{j-1},2^{-j}\}$.
The last inequality proves $\varepsilon_j\to0$.

For every compact set $K\Subset U$ and every $L$, the conditions with $N=0$, applied on a $K_p$ containing $K$ and with $p\geq L$, show that $\displaystyle \sum_{j=0}^{\infty}b_j$ converges in the seminorm $q^K_{L,m_0}$. Consequently, \[
a:=\sum_{j=0}^{\infty}b_j
\] belongs to $S^{m_0}$. Now fix $N$. The term $b_N$ belongs to $S^{m_N}$, and the preceding conditions, this time with target order $m_N$, show that $\displaystyle \sum_{j>N}b_j$ converges in all seminorms of $S^{m_N}$. Thus, for each $N$,
\[
a-\displaystyle\sum_{j=0}^{N-1}a_j
=
\displaystyle\sum_{j=N}^\infty b_j
+
\displaystyle\sum_{j=0}^{N-1}(b_j-a_j).
\]
The first sum belongs to $S^{m_N}$. Each difference $b_j-a_j$ has bounded support in $\xi$ and therefore belongs to $S^{-\infty}(U\times\mathbb R^n;
\operatorname{Hom}(\mathbb C^r,\mathbb C^s))$; this proves the asymptotic expansion. Ungrouping the finite groups recovers the interpretation specified in the definition. For the last clause, apply the same construction with radial cutoffs to $a_j=\vartheta h_{m_0-j}$; the remainder after the first $N$ terms belongs to $S^{m_0-N}$, which is exactly the expansion in Definition~\ref{def:simbolo-clasico-pseudo}.
\end{proof}

\subsection{Oscillatory integrals and quantization}

So far we have worked with symbol functions. To turn them into operators, we must give meaning to an integral that, except at very negative orders, is not absolutely convergent. Oscillation of the phase will compensate for growth of the symbol.

Let $n,r,s\in\mathbb N$, $m<-n$, $a\in S^m_{1,0}(\mathbb R^n;
\operatorname{Hom}(\mathbb C^r,\mathbb C^s))$, and $u\in\mathcal S(\mathbb R^n,\mathbb C^r)$. Then the integral
\[
(2\pi)^{-n}\int_{\mathbb R_\xi^n}\int_{\mathbb R_y^n}
e^{i(x-y)\cdot\xi}a(x,\xi)u(y)\,dy\,d\xi
\]
is absolutely convergent. At higher orders, oscillation of the exponential takes the place of absolute integrability. We make this replacement precise by regularization, rather than by an implicit convention.

Let $c\in C_c^\infty(\mathbb R^{2n})$ equal one near $(0,0)$. If $q(x,y,\xi)$ is smooth and all its derivatives have polynomial growth in $(y,\xi)$, define, whenever the limit exists,
\[
\operatorname{Os}\!\iint_{\mathbb R_y^n\times\mathbb R_\xi^n}
e^{i(x-y)\cdot\xi}q(x,y,\xi)\,dy\,d\xi
:=
\lim_{\varepsilon\to0^+}
\iint_{\mathbb R_y^n\times\mathbb R_\xi^n}
e^{i(x-y)\cdot\xi}c(\varepsilon y,\varepsilon\xi)
q(x,y,\xi)\,dy\,d\xi.
\]
Existence will be proved under growth bounds of fixed order. Specifically, suppose there exist $\mu,\tau\in\mathbb R$ such that, for each compact set $K$ in the variable $x$ and all multi-indices $\gamma,\alpha,\beta$,
\[
\|D_x^\gamma D_y^\alpha D_\xi^\beta q(x,y,\xi)\|
\leq C_{K,\gamma,\alpha,\beta}
\langle y\rangle^\tau\langle\xi\rangle^{\mu-|\beta|},
\qquad x\in K.
\]
The exponents $\mu$ and $\tau$ do not depend on the number of derivatives. The phase is the real function $(x-y)\cdot\xi$. The identities \[
e^{i(x-y)\cdot\xi}
=
\langle\xi\rangle^{-2N}(1-\Delta_y)^Ne^{i(x-y)\cdot\xi}
\] and \[
e^{i(x-y)\cdot\xi}
=
\langle x-y\rangle^{-2L}(1-\Delta_\xi)^Le^{i(x-y)\cdot\xi}
\] transfer derivatives to the amplitude. Write $c_\varepsilon(y,\xi)=c(\varepsilon y,\varepsilon\xi)$. Integrating by parts first in $y$ and then in $\xi$ shows that the regularized integral equals the integral of $e^{i(x-y)\cdot\xi}R_\varepsilon(x,y,\xi)$, where
\[
R_\varepsilon
=
\langle x-y\rangle^{-2L}(1-\Delta_\xi)^L
\bigl[\langle\xi\rangle^{-2N}(1-\Delta_y)^N
(c_\varepsilon q)\bigr].
\]
These operations are valid before taking the limit because the cutoff has compact support. For $0<\varepsilon\leq1$, its derivatives satisfy
\[
|D_y^\alpha D_\xi^\beta c_\varepsilon(y,\xi)|
\leq C_{\alpha,\beta}\langle\xi\rangle^{-|\beta|}.
\]
Indeed, each derivative contributes a factor of $\varepsilon$, and on the support of the corresponding term, $\varepsilon\langle\xi\rangle$ is bounded independently of $\varepsilon$. The Leibniz rule and the bounds on $q$ give, for $x\in K$,
\[
\|R_\varepsilon(x,y,\xi)\|
\leq C_K\langle y\rangle^{\tau-2L}
\langle\xi\rangle^{\mu-2N}.
\]
Here we use that $\langle x-y\rangle$ and $\langle y\rangle$ are uniformly comparable as $x$ ranges over $K$. Choose $2N>\mu+n$ and $2L>\tau+n$; the last expression is then integrable in $(y,\xi)$. For each fixed $(y,\xi)$, the cutoff equals one and all its positive-order derivatives vanish when $\varepsilon$ is sufficiently small. By dominated convergence, the limit is the absolutely convergent integral obtained by replacing $c_\varepsilon q$ with $q$ in $R_\varepsilon$. This expression is independent of the cutoff $c$.

Differentiating with respect to $x$ up to order $d$ also produces powers of $\xi$ of degree at most $d$. The same proof, with $2N>\mu+d+n$, gives an integrable majorant for all these derivatives, uniformly on $K$. Convergence of the integrand and its derivatives is also uniform in $x\in K$ for fixed $(y,\xi)$. This yields convergence in every $C^d(K)$ norm and smooth dependence of the oscillatory integral on $x$.

For the following quantization, $q(x,y,\xi)=a(x,\xi)u(y)$ satisfies these bounds with $\mu=m$ and $\tau=0$, since $u$ is a Schwartz function. The same holds for $b(x,y,\xi)u(y)$ for the fixed-order amplitudes introduced below. In the composition and reduction formulas, apply the argument to the integration variables, leaving the remaining variables as parameters; symbol estimates and Peetre's inequality control the weights appearing there. The bounds are local in the parameter variables. Polynomial growth of each derivative, without common control of the exponents, will not be used as a sufficient hypothesis for existence.

Together with asymptotic summation, regularization allows us to pass from formal expansions to well-defined operators. An equivalent construction, with remainders estimated seminorm by seminorm, can be found in \cite[Sections~1 and~3]{Shubin2001}.

\begin{definition}[Kohn--Nirenberg quantization]
\label{def:cuantizacion-kohn-nirenberg}
\index{quantization@quantization!Kohn--Nirenberg}
Let $n,r,s\in\mathbb N$, $m\in\mathbb R$, and $a\in S^m_{1,0}(\mathbb R^n;
\operatorname{Hom}(\mathbb C^r,\mathbb C^s))$. Define
\[
\operatorname{Op}(a)u(x)
:=
(2\pi)^{-n}
\operatorname{Os}\!\iint_{\mathbb R_y^n\times\mathbb R_\xi^n}
e^{i(x-y)\cdot\xi}a(x,\xi)u(y)\,dy\,d\xi,
\qquad u\in\mathcal S(\mathbb R^n,\mathbb C^r).
\]
Equivalently,
\[
\operatorname{Op}(a)u(x)
=
(2\pi)^{-\frac{n}{2}}
\operatorname{Os}\!\int_{\mathbb R_\xi^n}
e^{ix\cdot\xi}a(x,\xi)\widehat u(\xi)\,d\xi.
\]
\end{definition}

The second formula is first obtained when the integrals are absolutely convergent, using Fourier inversion and Fubini's theorem, and then by regularization. The factor $(2\pi)^{-n}$ in the first formula results from the two occurrences of the unitary normalization; it cannot be replaced by $(2\pi)^{-\displaystyle\frac{n}{2}}$.

\begin{proposition}[Continuity of quantization]
\label{prop:op-s-schwartz-continuidad}
Let $n,r,s\in\mathbb N$, $m\in\mathbb R$, and \[
a\in S^m_{1,0}(\mathbb R^n;
\operatorname{Hom}(\mathbb C^r,\mathbb C^s)).
\]. Then $\operatorname{Op}(a)$ is continuous from $\mathcal S(\mathbb R^n,\mathbb C^r)$ to $\mathcal S(\mathbb R^n,\mathbb C^s)$. By transposition it extends uniquely to a continuous operator
\[
\operatorname{Op}(a)\colon
\mathcal S'(\mathbb R^n,\mathbb C^r)
\longrightarrow
\mathcal S'(\mathbb R^n,\mathbb C^s).
\]
If $a\in S^{-\infty}(\mathbb R^n;
\operatorname{Hom}(\mathbb C^r,\mathbb C^s))$, it maps $\mathcal S'(\mathbb R^n,\mathbb C^r)$ into $C^\infty(\mathbb R^n,\mathbb C^s)$; if, in addition, $a$ decays rapidly also in $x$, it maps $\mathcal S'(\mathbb R^n,\mathbb C^r)$ into $\mathcal S(\mathbb R^n,\mathbb C^s)$.
\end{proposition}

\begin{proof}
Differentiating the second formula with respect to $x$ produces sums of terms with factors $\xi^\gamma D_x^{\alpha-\gamma}a(x,\xi)$. To control a power $x^\beta$, use $x_je^{ix\cdot\xi}=-iD_{\xi_j}e^{ix\cdot\xi}$ and integrate by parts in $\xi$. After expanding by the Leibniz rule, each Schwartz seminorm of $\operatorname{Op}(a)u$ is bounded by a finite sum of symbol seminorms of $a$ and Schwartz seminorms of $u$. Choose the number of integrations by parts greater than $m+n$ to make the variable $\xi$ integrable. This proves existence, membership in $\mathcal S(\mathbb R^n,\mathbb C^s)$, and continuity simultaneously.

The transposed operator is defined by $\langle\operatorname{Op}(a)u,\varphi\rangle
:=\langle u,\operatorname{Op}(a)^t\varphi\rangle$, where the transpose on the right-hand side is obtained by interchanging the kernel variables; the proposition on amplitudes proved below shows that it too preserves the corresponding Schwartz spaces. If $a\in S^{-\infty}(\mathbb R^n;
\operatorname{Hom}(\mathbb C^r,\mathbb C^s))$, the decay in $\xi$ can be chosen arbitrarily large, allowing differentiation under the integral against a tempered distribution and yielding smoothness. Additional decay in $x$ gives all Schwartz seminorms.
\end{proof}

\subsection{Amplitudes and reduction to symbols}

Kohn--Nirenberg quantization uses a symbol evaluated at the output variable $x$. Composing operators, taking adjoints, or changing coordinates naturally introduces coefficients that also depend on the input variable $y$. Before closing the calculus under these operations, we must show that such amplitudes can be reduced back to symbols.

\begin{definition}[Amplitude]
\label{def:amplitud-pseudodiferencial}
\index{amplitude!pseudodifferential}
Let $n,r,s\in\mathbb N$ and $m\in\mathbb R$. An \textbf{amplitude of order $m$} is a smooth map \[
b\colon\mathbb R_x^n\times\mathbb R_y^n\times\mathbb R_\xi^n
\longrightarrow\operatorname{Hom}(\mathbb C^r,\mathbb C^s)
\] such that, for all multi-indices $\alpha,\gamma,\beta$, there exists $C_{\alpha,\gamma,\beta}>0$ with
\[
\|D_x^\alpha D_y^\gamma D_\xi^\beta b(x,y,\xi)\|
\leq
C_{\alpha,\gamma,\beta}\langle\xi\rangle^{m-|\beta|}.
\]
Its operator, initially on $\mathcal S(\mathbb R^n,\mathbb C^r)$, is defined by
\[
I(b)u(x):=(2\pi)^{-n}
\operatorname{Os}\!\iint_{\mathbb R_y^n\times\mathbb R_\xi^n}
e^{i(x-y)\cdot\xi}b(x,y,\xi)u(y)\,dy\,d\xi.
\]

For an open set $U\subseteq\mathbb R^n$, a \textbf{local amplitude} on
$U\times U\times\mathbb R^n$ satisfies these estimates on every compact
$K\Subset U\times U$, with constants $C_{K,\alpha,\gamma,\beta}$.
It defines $I(b)$ initially on $C_c^\infty(U,\mathbb C^r)$.
We call it properly supported if its spatial support is contained in a
closed subset of $U\times U$ whose two projections to $U$ are proper.
In that case its operator acts on arbitrary distributions, by transposition.

An amplitude is \textbf{classical} if
$b\sim\sum_{j\geq0}b_{m-j}$, where $b_{m-j}$ is smooth for $\xi\neq0$,
homogeneous of degree $m-j$ in $\xi$, and the remainder after $N$ terms
(with a fixed cutoff near $\xi=0$) is an amplitude of order $m-N$.
The same definition applies with local estimates.
\end{definition}

A symbol is the special amplitude independent of $y$. The reverse direction contains the argument that will recur in composition, adjoints, and coordinate changes.

\begin{theorem}[Amplitude reduction]
\label{teo:amplitud-a-simbolo}
Let $n,r,s\in\mathbb N$ and $m\in\mathbb R$. For each amplitude $b$ of order $m$ with values in $\operatorname{Hom}(\mathbb C^r,\mathbb C^s)$, there exists \[
a\in S^m_{1,0}(\mathbb R^n;
\operatorname{Hom}(\mathbb C^r,\mathbb C^s))
\] such that $I(b)=\operatorname{Op}(a)$. The symbol can be chosen with expansion
\[
a(x,\xi)
\sim
\displaystyle\sum_{\alpha\in\mathbb N_0^n}
\frac{i^{-|\alpha|}}{\alpha!}
D_\xi^\alpha D_y^\alpha b(x,y,\xi)\big|_{y=x}.
\]
The exact global symbol is uniquely determined by the operator. For a
classical amplitude the reduced symbol is classical, with homogeneous terms
\[
 a_{m-k}(x,\xi)
 =\sum_{j+|\alpha|=k}\frac{i^{-|\alpha|}}{\alpha!}
   D_\xi^\alpha D_y^\alpha b_{m-j}(x,y,\xi)\big|_{y=x}.
\]
For local amplitudes the same expansion holds in
$S^m_{1,0}(U\times\mathbb R^n;\operatorname{Hom}(\mathbb C^r,\mathbb C^s))$,
modulo an operator with smooth kernel. Proper support gives the corresponding
identities on distributions.
\end{theorem}

\begin{proof}
Throughout this proof, $n,r,s,m$ remain fixed, and we write
\[
S^q:=S^q_{1,0}
\bigl(\mathbb R^n;
\operatorname{Hom}(\mathbb C^r,\mathbb C^s)\bigr).
\]
Exact reduction is obtained through the partial Fourier transform of the kernel in the variable $x-y$. With this chapter's normalization, define
\begin{equation}
\label{eq:simbolo-exacto-reduccion-amplitud}
a(x,\xi)
:=
(2\pi)^{-n}
\operatorname{Os}\!\iint_{\mathbb R_y^n\times\mathbb R_\eta^n}
e^{-iy\cdot\eta}
b(x,x+y,\xi+\eta)\,dy\,d\eta.
\end{equation}
This oscillatory integral exists. To see this, use, with $L,K\in\mathbb N$,
\[
e^{-iy\cdot\eta}
=
\langle\eta\rangle^{-2K}(1-\Delta_y)^K e^{-iy\cdot\eta},
\qquad
e^{-iy\cdot\eta}
=
\langle y\rangle^{-2L}(1-\Delta_\eta)^L e^{-iy\cdot\eta},
\]
applying the first identity first, and then the second to each resulting term. After integration by parts, the derivatives fall on the amplitude and the weight in $\eta$. Peetre's inequality in Lemma~\ref{lem:desigualdad-peetre-peso-japones} gives, for every $\mu\in\mathbb R$,
\[
\langle\xi+\eta\rangle^\mu
\leq 2^{\frac{|\mu|}{2}}\langle\xi\rangle^\mu
\langle\eta\rangle^{|\mu|}.
\]
Fix the multi-indices \(\alpha,\beta\) and choose integers
\(L,K\geq1\) with
\[
 2L>n,\qquad 2K>n+\bigl|m-|\beta|\bigr|+2L.
\]
After both integrations by parts, the factor multiplying
\(e^{-iy\cdot\eta}\) is
\[
 \mathcal R(x,y,\xi,\eta):=
 \langle y\rangle^{-2L}(1-\Delta_\eta)^L
 \bigl[\langle\eta\rangle^{-2K}(1-\Delta_y)^K
                         b(x,x+y,\xi+\eta)\bigr].
\]
An \(x\) derivative acts on both spatial variables of \(b\), whose
bounds are uniform. Each Leibniz term after applying
\(D_x^\alpha D_\xi^\beta\) contains at most \(2L\) additional
derivatives in the covariable of \(b\). If \(\gamma\) is their
multi-index, Peetre's inequality and \(\langle\xi\rangle\geq1\) give
\[
 \langle\xi+\eta\rangle^{m-|\beta|-|\gamma|}
 \leq C_{m,\beta,\gamma}\langle\xi\rangle^{m-|\beta|}
                  \langle\eta\rangle^{|m-|\beta||+2L}.
\]
Derivatives of the weight \(\langle\eta\rangle^{-2K}\) do not worsen
its decay. Summing the finitely many terms yields
\[
 \|D_x^\alpha D_\xi^\beta\mathcal R(x,y,\xi,\eta)\|
 \leq C_{\alpha,\beta,L,K}\langle\xi\rangle^{m-|\beta|}
       \langle y\rangle^{-2L}
       \langle\eta\rangle^{-2K+|m-|\beta||+2L}.
\]
The last two factors are integrable by the choice of \(L,K\).
The same majorant holds with the regularizing cutoff inserted:
its \(\eta\) derivatives contribute factors bounded by
\(C\langle\eta\rangle^{-d}\) at order \(d\), and its \(y\)
derivatives remain bounded. Dominated convergence allows removal of
the cutoff and differentiation under the integral. In particular,
\[
 \|D_x^\alpha D_\xi^\beta a(x,\xi)\|
 \leq C'_{\alpha,\beta}\langle\xi\rangle^{m-|\beta|}.
\]

The same argument for any number of derivatives proves that $a\in S^m$. Inserting compactly supported cutoffs before applying Fubini's theorem, the kernel formula and the changes $y-x\mapsto y$ and $\eta-\xi\mapsto\eta$ give $I(b)=\operatorname{Op}(a)$; the preceding estimates allow the cutoffs to be removed. Thus \eqref{eq:simbolo-exacto-reduccion-amplitud} directly constructs the symbol associated with the amplitude.

Now fix $N\geq1$ and expand the amplitude in $y$ about $x$:
\[
b(x,y,\xi)
=
\displaystyle\sum_{|\alpha|<N}
\frac{(y-x)^\alpha}{\alpha!}D_y^\alpha b(x,x,\xi)
+
\displaystyle\sum_{|\alpha|=N}(y-x)^\alpha r_\alpha(x,y,\xi),
\]
where
\[
r_\alpha(x,y,\xi)
=
\frac{N}{\alpha!}\int_0^1(1-t)^{N-1}
D_y^\alpha b(x,x+t(y-x),\xi)\,dt.
\]
Since $(y-x)^\alpha e^{i(x-y)\cdot\xi}=i^{|\alpha|}D_\xi^\alpha e^{i(x-y)\cdot\xi}$, integration by parts in $\xi$ turns each term of the finite sum into
\[
\frac{i^{-|\alpha|}}{\alpha!}
D_\xi^\alpha D_y^\alpha b(x,x,\xi).
\]
Equivalently, inserting this Taylor expansion into \eqref{eq:simbolo-exacto-reduccion-amplitud}, the identity \[
y^\alpha e^{-iy\cdot\eta}
=i^{|\alpha|}D_\eta^\alpha e^{-iy\cdot\eta}
\] produces the same terms. In the remainder, the $N$ integrations by parts make $D_\eta^\alpha$ act on $r_\alpha(x,x+y,\xi+\eta)$. Each derivative in $\eta$ is a derivative in the covariable of $b$ and lowers the order by one. The estimate used for \eqref{eq:simbolo-exacto-reduccion-amplitud}, now applied to these derivatives, shows directly that the reduced remainder belongs to $S^{m-N}$. Since $N$ is arbitrary, this gives the asserted expansion without a second appeal to amplitude reduction.

If $b$ is classical, fix $N$ and insert its expansion through degree
$m-N+1$ into the Taylor calculation. For $b_{m-j}$ keep exactly the terms
$|\alpha|<N-j$. All omitted terms have order at most $m-N$ by the remainder
estimate just proved. Grouping the retained terms by $k=j+|\alpha|$
gives the stated homogeneous coefficients; each has degree $m-k$.

In the global model, the partial Fourier transform of the kernel recovers
the exact symbol uniquely. For a local amplitude, choose a properly
supported cutoff equal to one near the diagonal. On a compact subset of
$U$, only a larger compact subset of $U\times U$ is then relevant, so the
same estimates apply after spatial cutoffs and extension to the global
model. The discarded kernel is smooth by integration by parts away from
the diagonal. These local reductions agree modulo smoothing symbols on
overlaps, giving the asserted local class and expansion. Choices of the
diagonal cutoff change the symbol by
$S^{-\infty}(U\times\mathbb R^n;\operatorname{Hom}(\mathbb C^r,\mathbb C^s))$;
this is a local assertion and requires no uniform bound as $x$ leaves
compact subsets of $U$.
\end{proof}

\subsection{Fundamental examples}

Before developing the algebraic calculus, it is useful to recognize some familiar operators within the new class. The following examples show how symbol order records loss or gain of derivatives while distinguishing locality, nonlocality, and smoothing.

\begin{example}[Differential operators]
Let $U\subseteq\mathbb R^n$ be open, and let $r,s\in\mathbb N$ and $A_\alpha\in C^\infty(U;
\operatorname{Hom}(\mathbb C^r,\mathbb C^s))$. For \[
P=\displaystyle\sum_{|\alpha|\leq m}A_\alpha(x)D^\alpha
\], the Fourier symbol is
\[
p(x,\xi)=\displaystyle\sum_{|\alpha|\leq m}A_\alpha(x)(i\xi)^\alpha,
\qquad
P=\operatorname{Op}(p).
\]
If $U=\mathbb R^n$ and the coefficients and all their derivatives are bounded, \[
p\in S^m_{1,0}(\mathbb R^n;
\operatorname{Hom}(\mathbb C^r,\mathbb C^s));
\]; locally it always belongs to $S^m_{1,0}(U;\operatorname{Hom}(\mathbb C^r,\mathbb C^s))$. The principal Fourier component is
\[
p_m(x,\xi)=i^m\displaystyle\sum_{|\alpha|=m}A_\alpha(x)\xi^\alpha
=i^m\boldsymbol{\sigma}_m(P)(x,\xi),
\]
where $\boldsymbol{\sigma}_m(P)$ is exactly the intrinsic symbol of Definition~\ref{def:simbolo-principal-pdo-haces}. This reconciles the previously fixed differential convention with Fourier quantization. Multiplication by $i^m$ does not change ellipticity.

The operator is local: its kernel is a finite sum of transverse derivatives of the delta distribution on the diagonal,
\[
\mathbf{K}_P(x,y)=\sum_{|\alpha|\leq m}A_\alpha(x)D_x^\alpha\delta(x-y).
\]
Its support is therefore contained in the diagonal, where all its singularities are concentrated. The Fourier and localization estimates already proved in the Sobolev chapters give $P\colon H^{t+m}_{\mathrm{loc}}(U,\mathbb C^r)\longrightarrow
H^t_{\mathrm{loc}}(U,\mathbb C^s)$; if $U=\mathbb R^n$ and the coefficients and their derivatives are globally bounded, one obtains $P\colon H^{t+m}(\mathbb R^n,\mathbb C^r)\longrightarrow
H^t(\mathbb R^n,\mathbb C^s)$.
\end{example}

\begin{example}[Multipliers and Bessel potentials]
If $a$ depends only on $\xi$, then $\operatorname{Op}(a)u=\mathcal F^{-1}(a\widehat u)$. For $s\in\mathbb R$, the estimates preceding Definition~\ref{def:regularidad-intermedia-norma-espacio-sobolev-fraccionario} show that $\langle\xi\rangle^s\in
S^s_{\mathrm{cl}}(\mathbb R^n)$, with principal component $\|\xi\|^s$ for $\xi\neq0$. Thus
\[
J^s=\operatorname{Op}(\langle\xi\rangle^s),
\qquad
J^{-s}=\operatorname{Op}(\langle\xi\rangle^{-s}).
\]
The first has order $s$ and the second order $-s$; neither is local unless its symbol is a polynomial. When $s>0$, the kernel of $J^{-s}$ is the Bessel kernel in Proposition~\ref{prop:nucleo-bessel-euclidiano}.

Their principal symbols are $\|\xi\|^s$ and $\|\xi\|^{-s}$, respectively. Directly from the Fourier definition, for every $t\in\mathbb R$,
\[
J^s\colon H^{t+s}(\mathbb R^n)\longrightarrow H^t(\mathbb R^n),
\qquad
J^{-s}\colon H^t(\mathbb R^n)\longrightarrow
H^{t+s}(\mathbb R^n)\quad(s>0)
\]
are isomorphisms. The kernels are smooth when $x\neq y$ and, except in the local polynomial cases, have a conormal singularity at $x=y$; for $s>0$, the kernel of $J^{-s}$ progressively loses this singularity as $s$ increases, but does not become $C^\infty$ at a general finite order.
\end{example}

\begin{example}[The fractional Laplacian and low frequencies]
Let $s>0$. The multiplier of $(-\Delta)^{\displaystyle\frac{s}{2}}$ is $\|\xi\|^s$. If $s$ is not an even integer, it is not smooth at $\xi=0$ and therefore does not belong to the global inhomogeneous class $S^s_{1,0}(\mathbb R^n)$ just defined. If $\vartheta$ vanishes near zero and equals one for $\|\xi\|\geq1$, then $\vartheta(\xi)\|\xi\|^s\in S^s_{\mathrm{cl}}(\mathbb R^n)$. The remaining part has compact frequency support and maps compactly supported distributions to smooth functions. Thus the fractional Laplacian is microlocally a classical pseudodifferential operator of order $s$, although its global homogeneous treatment retains additional information at frequency zero. This distinction disappears for $(I-\Delta)^{\displaystyle\frac{s}{2}}$, whose symbol $\langle\xi\rangle^s$ is smooth on all of $\mathbb R^n$.

For $s>0$, the multiplier $\|\xi\|^s$ maps $H^{t+s}(\mathbb R^n)$ continuously into $H^t(\mathbb R^n)$ and has principal symbol $\|\xi\|^s$. It is nonlocal when $s$ is not an even integer. If $0<s<2$, its kernel off the diagonal is represented, up to a dimensional constant, by $\|x-y\|^{-n-s}$, and its action is interpreted in the principal-value sense through the difference $u(x)-u(y)$. Thus the kernel's singularity lies on the diagonal, but its support does not: this is the distinction between nonlocality and pseudolocality.
\end{example}

\begin{example}[Riesz and Hilbert transforms]
The Riesz transform $R_j=D_j(-\Delta)^{-\displaystyle\frac{1}{2}}$ has symbol $\displaystyle\frac{i\xi_j}{\|\xi\|}$ away from $\xi=0$, homogeneous of degree zero. In dimension one, the Hilbert transform with principal-value kernel $\pi^{-1}(x-y)^{-1}$ has multiplier $-i\operatorname{sign}(\xi)$. Both are nonlocal operators of order zero at high frequencies; a cutoff at frequency zero produces representatives in $S^0_{\mathrm{cl}}(\mathbb R^n)$ with the same microlocal behavior.

These multipliers are also their principal symbols. Since they are bounded, both operators act continuously on every $H^t(\mathbb R^n)$. Away from the diagonal, the Riesz kernel is a multiple of $\displaystyle\frac{x_j-y_j}{\|x-y\|^{n+1}}$, and the Hilbert kernel is $\pi^{-1}\operatorname{pv}(x-y)^{-1}$. They are smooth away from $x=y$, singular on the diagonal, and not supported there.
\end{example}

\begin{example}[Smoothing operators and exact inverses]
The symbol $e^{-\|\xi\|^2}$ belongs to $S^{-\infty}(\mathbb R^n;
\operatorname{Hom}(\mathbb C,\mathbb C))$, and its operator is convolution with a Gaussian; it maps $\mathcal S'(\mathbb R^n)$ into $C^\infty(\mathbb R^n)$. On the other hand, $I-\Delta$ has multiplier $1+\|\xi\|^2$, and its exact inverse on $\mathcal S'(\mathbb R^n)$ is $\operatorname{Op}(\langle\xi\rangle^{-2})$. This last example contains the general mechanism of elliptic symbolic inversion without error terms.

The Gaussian operator has order $-\infty$, zero principal symbol in every finite-order quotient, and a smooth nonlocal kernel, and it maps $H^t(\mathbb R^n)$ into $H^N(\mathbb R^n)$ for all $t,N\in\mathbb R$. The inverse $(I-\Delta)^{-1}$ has order $-2$ and principal symbol $\|\xi\|^{-2}$ and is an isomorphism $H^t(\mathbb R^n)\longrightarrow H^{t+2}(\mathbb R^n)$. Its Bessel kernel is smooth off the diagonal and has the singularity there corresponding to order $-2$; it is an exact inverse and, in particular, a parametrix without a remainder.
\end{example}

\section{Kernels, proper support, and symbolic calculus}
\label{sec:pseudo-nucleos-calculo}

The symbol describes the operator in frequency, whereas the kernel shows how two points in space interact. This second description is essential for distinguishing local and nonlocal operators, formulating proper support, and proving pseudolocality. It also explains why composition, adjoints, and coordinate changes remain within the calculus.

\subsection{The kernel and the diagonal}

The Schwartz kernel theorem for vector bundles, Theorem~\ref{teo: nucleo Schwartz haces}, represents a continuous operator by a distribution on the product. In the Euclidean scalar model, Definition~\ref{def:cuantizacion-kohn-nirenberg} gives
\[
\mathbf{K}_{\operatorname{Op}(a)}(x,y)
=
(2\pi)^{-n}\operatorname{Os}\!\int_{\mathbb R^n}
e^{i(x-y)\cdot\xi}a(x,\xi)\,d\xi.
\]
The possible singularity is concentrated at $x=y$. If $x\neq y$, the operator \[
L_{x,y}
:=
\frac{1}{i\|x-y\|^2}(x-y)\cdot D_\xi
\] satisfies $L_{x,y}e^{i(x-y)\cdot\xi}=e^{i(x-y)\cdot\xi}$. After integrating by parts $N$ times, derivatives of order $N$ with respect to $\xi$ reduce the symbol to order $m-N$. For $N>m+n$, the resulting integral is absolutely convergent. The same procedure may be repeated after differentiating in $x$ and $y$. Thus $\mathbf{K}_{\operatorname{Op}(a)}$ is smooth off the diagonal
\[
\Delta_{\mathbb R^n}:=\{(x,x)\mid x\in\mathbb R^n\}.
\]

An operator $A$ is called \textbf{pseudolocal} if, for every distribution $u$ on which it is defined,
\[
\operatorname{sing\,supp}(Au)
\subseteq
\operatorname{sing\,supp}(u).
\]
This property does not say that $A$ is local: it only prevents it from creating singularities at points where the datum is smooth.

Thus a pseudodifferential operator can be nonlocal because its kernel is not supported on the diagonal, while remaining pseudolocal because that kernel is smooth off the diagonal. Extending the operator to general distributions and controlling supports also requires a global condition on the projections of its kernel.

\begin{definition}[Proper support]
\label{def:soporte-propio-pseudo}
\index{proper support!pseudodifferential operator}
Let $N$ and $M$ be smooth manifolds without boundary, let $\mathbf{E}\to N$ and $\mathbf{F}\to M$ be smooth vector bundles of finite rank, let \[
A\colon \Gamma_c(N,\mathbf{E})\longrightarrow\Gamma(M,\mathbf{F})
\] be a continuous operator, also regarded as taking values in $\mathcal D'(M,\mathbf{F})$, and let $\mathbf{K}_A$ be its kernel on $M\times N$. We say that $A$ is \textbf{properly supported} if both projections
\[
\operatorname{pr}_M\colon\supp \mathbf{K}_A\longrightarrow M,
\qquad
\operatorname{pr}_N\colon\supp \mathbf{K}_A\longrightarrow N
\]
are proper maps. Equivalently, for every compact subset $K$ of either factor, the intersection of $\supp \mathbf{K}_A$ with the preimage of $K$ has compact projection onto the other factor.
\end{definition}

In this chapter's pseudodifferential calculus we take $N=M$; the formulation with two manifolds simply records that proper support is a condition on the kernel and its two projections.

This condition has two distinct consequences. Projection onto the source, $\operatorname{pr}_N$, controls the support of $Au$ when $u$ has compact support. Projection onto the target, $\operatorname{pr}_M$, allows $Au$ to be defined for a distribution without compact support: local evaluation against a test function involves only a compact part of the distribution. In particular, a properly supported operator extends continuously
\[
A\colon \Gamma_c(\mathbf{E})\longrightarrow \Gamma_c(\mathbf{F}),
\qquad
A\colon\mathcal D'(N,\mathbf{E})\longrightarrow\mathcal D'(M,\mathbf{F}),
\]
if its local expression is pseudodifferential. The second extension is defined by duality with the properly supported transpose.

\begin{proposition}[Decomposition with proper support]
\label{prop:descomposicion-soporte-propio-pseudo}
Let $U\subseteq\mathbb R^n$ be open, and let $r,s\in\mathbb N$ and $a\in S^m_{1,0}(U;
\operatorname{Hom}(\mathbb C^r,\mathbb C^s))$. For the local operator $A=\operatorname{Op}(a)$, there exist a properly supported pseudodifferential operator $A_0$ of the same order and an operator $R$ with smooth kernel such that $A=A_0+R$. Moreover, $A_0$ has the same principal symbol as $A$.
\end{proposition}

\begin{proof}
Choose a function $\chi(x,y)$ equal to one on a neighborhood of the diagonal and supported in a properly supported neighborhood of it. Define $\mathbf{K}_{A_0}:=\chi \mathbf{K}_A$. Then $A_0$ is properly supported. The amplitude $\chi(x,y)a(x,\xi)$ reduces to a symbol by Theorem~\ref{teo:amplitud-a-simbolo}. In its expansion, the term with $\alpha=0$ is $a(x,\xi)$, and for $|\alpha|\geq1$,
\[
D_\xi^\alpha D_y^\alpha
\bigl(\chi(x,y)a(x,\xi)\bigr)\big|_{y=x}=0,
\]
since $\chi$ is identically one on a neighborhood of the diagonal. Thus the reduced symbol differs from $a$ by an element of $S^{-\infty}(U\times\mathbb R^n;
\operatorname{Hom}(\mathbb C^r,\mathbb C^s))$ and, in particular, has the same principal symbol. On the other hand, $(1-\chi)\mathbf{K}_A$ is supported outside a neighborhood of the diagonal. The integrations by parts at the beginning of this section show that it is smooth. Thus $R=A-A_0$ has the asserted kernel.
\end{proof}

\begin{definition}[Smoothing operator]
\label{def:operador-regularizante-pseudo}
\index{smoothing operator}
Let $U\subseteq\mathbb R^n$ be open and let $r,s\in\mathbb N$. The class \[
\Psi^{-\infty}_{\mathrm{loc}}
(U;\mathbb C^r,\mathbb C^s)
\] consists of continuous linear operators \[
A\colon C_c^\infty(U,\mathbb C^r)\longrightarrow
C^\infty(U,\mathbb C^s)
\] whose Schwartz kernel, written relative to Lebesgue measure in the input variable, belongs to
\[
C^\infty\!\left(
U\times U;\operatorname{Hom}(\mathbb C^r,\mathbb C^s)
\right).
\]
Its elements are called \textbf{local smoothing operators}. Define
\begin{equation}
\Psi^{-\infty}(U;\mathbb C^r,\mathbb C^s)
:=
\left\{
A\in\Psi^{-\infty}_{\mathrm{loc}}(U;\mathbb C^r,\mathbb C^s)
\;\middle|\;
\mathbf{K}_A\text{ is properly supported}
\right\}.
\label{eq:def-operadores-regularizantes-euclidianos-propios}
\end{equation}
Thus the notation without the subscript $\mathrm{loc}$ includes proper support.

On all of $\mathbb R^n$, if \[
a\in S^{-\infty}(\mathbb R^n;
\operatorname{Hom}(\mathbb C^r,\mathbb C^s)),
\], the global Kohn--Nirenberg operator $\operatorname{Op}(a)$ has smooth kernel and, by Proposition~\ref{prop:op-s-schwartz-continuidad}, acts continuously from $\mathcal S'(\mathbb R^n,\mathbb C^r)$ to $C^\infty(\mathbb R^n,\mathbb C^s)$. This global statement imposes no proper-support condition; the operator belongs to $\Psi^{-\infty}(\mathbb R^n;\mathbb C^r,\mathbb C^s)$ only if its kernel also satisfies that condition.
\end{definition}

\begin{proposition}[Smoothing symbols and kernels]
\label{prop:equiv-simbolo-nucleo-regularizante}
Let $U\subseteq\mathbb R^n$ be open, let $r,s\in\mathbb N$, and let $A$ be a local pseudodifferential operator from $C_c^\infty(U,\mathbb C^r)$ to $C^\infty(U,\mathbb C^s)$. The following conditions are equivalent:
\begin{enumerate}[label=(\alph*)]
\item every localization of $A$ has a symbol in $S^{-\infty}(U\times\mathbb R^n;
\operatorname{Hom}(\mathbb C^r,\mathbb C^s))$;
\item $\mathbf{K}_A$ is smooth near the diagonal;
\item $\mathbf{K}_A$ belongs to $C^\infty(U\times U;
\operatorname{Hom}(\mathbb C^r,\mathbb C^s))$, that is, $A\in\Psi^{-\infty}_{\mathrm{loc}}(U;\mathbb C^r,\mathbb C^s)$.
\end{enumerate}
If $A$ is properly supported, these conditions are also equivalent to $A\in\Psi^{-\infty}(U;\mathbb C^r,\mathbb C^s)$.
\end{proposition}

\begin{proof}
If a localization has symbol in $S^{-\infty}(U\times\mathbb R^n;
\operatorname{Hom}(\mathbb C^r,\mathbb C^s))$, applying a total of $q$ derivatives in $x$ and $y$ to the kernel produces only polynomial factors in $\xi$ of degree at most $q$. Estimate \eqref{eq:estimaciones-simbolo-regularizante-euclidiano} with $N>n+q$ yields an integrable majorant in $\mathbb R^n_\xi$; differentiation under the integral is therefore allowed, and the kernel is smooth near the diagonal. Every pseudodifferential kernel is smooth off the diagonal, as proved at the beginning of this section, so (b) implies (c).

Conversely, suppose (c) holds and choose $\chi(x,y)$ equal to one on a neighborhood of the diagonal, with properly supported support inside a coordinate localization. The partial transform
\[
a_\chi(x,\xi)
:=
\int_U e^{-i(x-y)\cdot\xi}\chi(x,y)\mathbf{K}_A(x,y)\,dy
\]
is a rapidly decaying symbol. For multi-indices $\alpha,\beta$, first use $(D_x+D_y)e^{-i(x-y)\cdot\xi}=0$ and integration by parts in $y$ to obtain
\[
D_x^\alpha D_\xi^\beta a_\chi(x,\xi)
=(-i)^{|\beta|}\int_U e^{-i(x-y)\cdot\xi}
(D_x+D_y)^\alpha
\bigl((x-y)^\beta\chi(x,y)\mathbf{K}_A(x,y)\bigr)\,dy.
\]
Since $(1-\Delta_y)^Le^{-i(x-y)\cdot\xi}
=\langle\xi\rangle^{2L}e^{-i(x-y)\cdot\xi}$, a second integration by parts gives, for each compact set $K$ in the variable $x$,
\[
\sup_{x\in K}
\|D_x^\alpha D_\xi^\beta a_\chi(x,\xi)\|
\leq C_{K,\alpha,\beta,L}\langle\xi\rangle^{-2L}.
\]
Proper support of $\chi$ ensures that the variable $y$ ranges over a compact set when $x\in K$. Since $L$ is arbitrary, $a_\chi\in S^{-\infty}(U\times\mathbb R^n;
\operatorname{Hom}(\mathbb C^r,\mathbb C^s))$. The corresponding operator has kernel $\chi \mathbf{K}_A$. The part $(1-\chi)\mathbf{K}_A$ is smooth for every pseudodifferential operator by the off-diagonal smoothness proved above. Thus every localization has a smoothing symbol, giving (a). The last assertion follows directly from Definition~\ref{def:operador-regularizante-pseudo} upon adding proper support.
\end{proof}

\begin{theorem}[Pseudolocality]
\label{teo:pseudolocalidad-pseudo}
\index{pseudolocal operator}
If $A$ is pseudodifferential and properly supported, then
\[
\operatorname{sing\,supp}(Au)
\subseteq
\operatorname{sing\,supp}(u)
\]
for every distribution $u$.
\end{theorem}

\begin{proof}
Let $x_0\notin\operatorname{sing\,supp}(u)$. Choose $\chi\in C_c^\infty$ equal to one on a neighborhood of the closure of an open set $W\ni x_0$ and supported in a region where $u$ is smooth. Write $u=\chi u+(1-\chi)u$. The first term is a smooth compactly supported function, and $A(\chi u)$ is smooth by continuity of quantization. If $x\in W$ and $y\in\operatorname{supp}(1-\chi)$, then $x\neq y$; hence the pairs involved in $A((1-\chi)u)$ remain off the diagonal, where the kernel is smooth. Proper support allows $y$ to be restricted to a compact set when evaluating locally, so pairing this smooth kernel with $(1-\chi)u$ depends smoothly on $x$. Thus $Au$ is smooth on $W$.
\end{proof}

\subsection{Composition}

Amplitude reduction now closes the calculus under composition. The ordinary product of the symbols gives the leading term; mixed derivatives in $x$ and $\xi$ describe lower-order corrections.

Let $r_0,r_1,r_2\in\mathbb N$, $A=\operatorname{Op}(a)$, and $B=\operatorname{Op}(b)$, with
\[
a\in S^m_{1,0}(\mathbb R^n;
\operatorname{Hom}(\mathbb C^{r_1},\mathbb C^{r_2})),
\qquad
b\in S^\ell_{1,0}(\mathbb R^n;
\operatorname{Hom}(\mathbb C^{r_0},\mathbb C^{r_1})).
\]
As operators on the Schwartz class, $A$ and $B$ can be composed without additional hypotheses. Substituting the integral formulas, first inserting compactly supported cutoffs and regrouping the phases, yields an oscillatory integral with an amplitude. After the change $z=x+y$ and $\eta=\zeta-\xi$, the composition symbol can be written as
\[
(a\# b)(x,\xi)
=
(2\pi)^{-n}
\operatorname{Os}\!\iint_{\mathbb R_y^n\times\mathbb R_\eta^n}
e^{-iy\cdot\eta}a(x,\xi+\eta)b(x+y,\xi)\,dy\,d\eta.
\]
The Taylor expansion of $b(x+y,\xi)$ in $y=0$, followed by integration by parts in $\eta$, gives the following result.

\begin{theorem}[Symbolic composition]
\label{teo:composicion-simbolica-pseudo}
\index{symbolic calculus@symbolic calculus!composition}
Let $n,r_0,r_1,r_2\in\mathbb N$ and $m,\ell\in\mathbb R$, and let
\[
a\in S^m_{1,0}(\mathbb R^n;
\operatorname{Hom}(\mathbb C^{r_1},\mathbb C^{r_2})),
\qquad
b\in S^\ell_{1,0}(\mathbb R^n;
\operatorname{Hom}(\mathbb C^{r_0},\mathbb C^{r_1})).
\]
As operators from $\mathcal S(\mathbb R^n,\mathbb C^{r_0})$ to $\mathcal S(\mathbb R^n,\mathbb C^{r_2})$,
\[
\operatorname{Op}(a)\operatorname{Op}(b)
=
\operatorname{Op}(a\# b),
\qquad
a\# b\in S^{m+\ell}_{1,0}(\mathbb R^n;
\operatorname{Hom}(\mathbb C^{r_0},\mathbb C^{r_2})),
\]
and
\[
a\# b
\sim
\displaystyle\sum_{\alpha\in\mathbb N_0^n}
\frac{i^{-|\alpha|}}{\alpha!}
(D_\xi^\alpha a)(D_x^\alpha b).
\]
In particular, $(a\# b)_{m+\ell}=a_mb_\ell$ for classical symbols. The same formulas apply to operators localized on an open set and to their extensions to distributions whenever proper support ensures that the composition is defined.
\end{theorem}

\begin{proof}
The preceding integral formula is first obtained with compactly supported cutoffs, where Fubini's theorem applies, and the cutoffs are then removed using oscillatory integral estimates. For $N\geq1$, write
\[
b(x+y,\xi)
=
\displaystyle\sum_{|\alpha|<N}
\frac{y^\alpha}{\alpha!}D_x^\alpha b(x,\xi)
+
\displaystyle\sum_{|\alpha|=N}y^\alpha r_\alpha(x,y,\xi).
\]
Here we use the integral Taylor remainder
\[
r_\alpha(x,y,\xi)
=
\frac{N}{\alpha!}\int_0^1(1-t)^{N-1}
D_x^\alpha b(x+ty,\xi)\,dt,
\qquad |\alpha|=N.
\]
Since $y^\alpha e^{-iy\cdot\eta}=i^{|\alpha|}D_\eta^\alpha e^{-iy\cdot\eta}$, integration by parts makes $i^{-|\alpha|}D_\eta^\alpha$ act on $a(x,\xi+\eta)$. Integrating $e^{-iy\cdot\eta}$ against a function independent of $y$ recovers evaluation at $\eta=0$, producing the displayed terms.

The remaining part is a finite sum, for $|\alpha|=N$, of integrals in $t\in[0,1]$ whose amplitudes are
\[
\frac{N i^{-N}}{\alpha!}(1-t)^{N-1}
(D_\xi^\alpha a)(x,\xi+\eta)
(D_x^\alpha b)(x+ty,\xi).
\]
These amplitudes have order $m+\ell-N$ uniformly in $t$. The reduction estimate proved in Theorem~\ref{teo:amplitud-a-simbolo} provides, for each seminorm of order $P$, another of finite order $P'$ and a constant independent of $t$ such that
\[
q_{P,m+\ell-N}(r_N^\#)
\leq
C_{P,N,m,\ell}\,q_{P',m}(a)q_{P',\ell}(b).
\]
On an open set, use the same inequality on $K\Subset U$, replacing $K$ on the right-hand side by a slightly larger compact set. Thus, in the local case, $r_N^\#\in S^{m+\ell-N}_{1,0}
(U\times\mathbb R^n;
\operatorname{Hom}(\mathbb C^{r_0},\mathbb C^{r_2}))$. Since $N$ is arbitrary, Lemma~\ref{lem:suma-asintotica-simbolos} completes the expansion.
\end{proof}

\begin{corollary}[Commutators]
\label{cor:simbolo-conmutador-pseudo}
Let $a\in S^m_{1,0}(\mathbb R^n)$ and $b\in S^\ell_{1,0}(\mathbb R^n)$ be scalar symbols, and suppose composition is defined in both orders. Then $[\operatorname{Op}(a),\operatorname{Op}(b)]$ has order at most $m+\ell-1$. More precisely,
\[
a\#b-b\#a=\frac{1}{i}\{a,b\}+r,
\qquad r\in S^{m+\ell-2}_{1,0}(\mathbb R^n),
\]
where
\[
\{a,b\}:=\sum_{j=1}^n
\left(\frac{\partial a}{\partial\xi_j}\frac{\partial b}{\partial x_j}
-\frac{\partial a}{\partial x_j}\frac{\partial b}{\partial\xi_j}\right).
\]
If, in addition, $a$ and $b$ are classical, the commutator is classical and its homogeneous component of order $m+\ell-1$ is
\[
\frac{1}{i}\{a_m,b_\ell\}
:=
\frac{1}{i}\displaystyle\sum_{j=1}^n
\left(
\frac{\partial a_m}{\partial\xi_j}\frac{\partial b_\ell}{\partial x_j}
-
\frac{\partial a_m}{\partial x_j}\frac{\partial b_\ell}{\partial\xi_j}
\right).
\]
If \(r\in\mathbb N\),
\[
a\in S^m_{\mathrm{cl}}\bigl(\mathbb R^n;M_{r\times r}(\mathbb C)\bigr),
\qquad
b\in S^\ell_{\mathrm{cl}}\bigl(\mathbb R^n;M_{r\times r}(\mathbb C)\bigr),
\]
then the term of order \(m+\ell\) is \([a_m,b_\ell]\). If it vanishes, the commutator has order at most \(m+\ell-1\), and its homogeneous component of that degree is
\[
[a_m,b_{\ell-1}]
+
[a_{m-1},b_\ell]
+
\frac{1}{i}\sum_{j=1}^n
\left(
\frac{\partial a_m}{\partial\xi_j}
\frac{\partial b_\ell}{\partial x_j}
-
\frac{\partial b_\ell}{\partial\xi_j}
\frac{\partial a_m}{\partial x_j}
\right).
\]
In particular, if \(a_m\) and \(b_\ell\) are scalar multiples of the identity, the two algebraic commutators in this formula vanish, leaving only the Poisson term.
\end{corollary}
\begin{proof}
The composition theorem gives
\[
a\# b
\sim
\sum_{\alpha\in\mathbb N_0^n}
\frac{i^{-|\alpha|}}{\alpha!}
(D_\xi^\alpha a)(D_x^\alpha b),
\qquad
b\# a
\sim
\sum_{\alpha\in\mathbb N_0^n}
\frac{i^{-|\alpha|}}{\alpha!}
(D_\xi^\alpha b)(D_x^\alpha a).
\]
In the scalar case, take the terms with $|\alpha|<2$ in both expansions. The remainders belong to $S^{m+\ell-2}_{1,0}$, and the products without derivatives cancel because $ab=ba$. Thus
\[
a\#b-b\#a
=\frac{1}{i}\sum_{j=1}^n
\bigl[(D_{\xi_j}a)(D_{x_j}b)-(D_{\xi_j}b)(D_{x_j}a)\bigr]+r
=\frac{1}{i}\{a,b\}+r,
\]
with $r\in S^{m+\ell-2}_{1,0}$. Each term in the sum has order at most $m+\ell-1$, proving the assertion without assuming classicality.

If the scalar symbols are classical, at high frequencies we may write $a=a_m+a_{<m}$ and $b=b_\ell+b_{<\ell}$, where the remainders have orders $m-1$ and $\ell-1$. Poisson terms containing either remainder have order at most $m+\ell-2$. Thus the homogeneous component of order $m+\ell-1$ is $\displaystyle\frac{1}{i}\{a_m,b_\ell\}$. Composition of the classical expansions likewise gives the homogeneous components of lower orders.

For matrix-valued symbols, the same grouping must preserve the order of the factors. Degree \(m+\ell\) is \(a_mb_\ell-b_\ell a_m=[a_m,b_\ell]\). If it vanishes, the next degree is the sum of the derivative-free terms \[
[a_m,b_{\ell-1}]+[a_{m-1},b_\ell]
\] and the terms with one derivative displayed in the statement. If both principal symbols are scalar multiples of the identity, these algebraic commutators vanish.
\end{proof}

\subsection{Transpose and formal adjoint}

Composition shows that the calculus is stable under chaining operators. The second operation needed for elliptic theory is the adjoint. At the kernel level, transposition interchanges $x$ and $y$; in the Hermitian case one also takes the fiberwise adjoint. The resulting amplitude depends on $y$ and is reduced by Theorem~\ref{teo:amplitud-a-simbolo}. We denote by $\operatorname{Op}(a)_h^*$ the formal Hermitian adjoint determined by the pairing of Schwartz sections.

\begin{theorem}[Symbol of the adjoint]
\label{teo:simbolo-adjunto-pseudo}
\index{symbolic calculus@symbolic calculus!adjoint}
Let $n,r,s\in\mathbb N$, $m\in\mathbb R$, and \[
a\in S^m_{1,0}(\mathbb R^n;
\operatorname{Hom}(\mathbb C^r,\mathbb C^s)).
\]. With respect to the inner products on $L^2(\mathbb R^n,\mathbb C^r)$ and $L^2(\mathbb R^n,\mathbb C^s)$ defined by Lebesgue measure and the standard Hermitian products,
\[
\operatorname{Op}(a)_h^*=\operatorname{Op}(a^{(*)}),
\]
where \[
a^{(*)}\in S^m_{1,0}(\mathbb R^n;
\operatorname{Hom}(\mathbb C^s,\mathbb C^r))
\] and
\[
a^{(*)}(x,\xi)
\sim
\displaystyle\sum_{\alpha\in\mathbb N_0^n}
\frac{i^{-|\alpha|}}{\alpha!}
D_\xi^\alpha D_x^\alpha\bigl(a(x,\xi)^*\bigr).
\]
In particular,
\[
a^{(*)}-a^*\in S^{m-1}_{1,0}(\mathbb R^n;
\operatorname{Hom}(\mathbb C^s,\mathbb C^r)).
\]
If $a$ is classical, then so is $a^{(*)}$, and $(a^{(*)})_m=a_m^*$.
\end{theorem}

\begin{proof}
For $u\in\mathcal S(\mathbb R^n,\mathbb C^r)$ and $v\in\mathcal S(\mathbb R^n,\mathbb C^s)$, the kernel formula and regularized Fubini theorem give
\[
\langle\operatorname{Op}(a)u,v\rangle
=
(2\pi)^{-n}
\operatorname{Os}\!\iiint_{\mathbb R_x^n\times\mathbb R_y^n\times\mathbb R_\xi^n}
\left\langle u(y),e^{i(y-x)\cdot\xi}a(x,\xi)^*v(x)\right\rangle
\,dx\,dy\,d\xi.
\]
Thus the adjoint amplitude is $a(y,\xi)^*$ after renaming the variables. Its reduction to a symbol is exactly the asserted expansion. The term with $\alpha=0$ is $a^*$, and the remainder has order at most $m-1$. If $a$ is classical, substituting its homogeneous expansion and grouping terms of equal degree yields a classical expansion for $a^{(*)}$: each derivative in $\xi$ lowers the degree by one, and at degree $m$ only $a_m^*$ remains.
\end{proof}

For a differential operator of integer order $m$, the principal Fourier component is $a_m=i^m\boldsymbol{\sigma}_m(P)$. Thus
\[
(a_m)^*=(-i)^m\boldsymbol{\sigma}_m(P)^*
=i^m\bigl((-1)^m\boldsymbol{\sigma}_m(P)^*\bigr).
\]
Removing the conversion factor $i^m$ gives exactly \[
\boldsymbol{\sigma}_m(P_h^*)=(-1)^m\boldsymbol{\sigma}_m(P)^*
\] in Proposition~\ref{prop:simbolo-adjunto-formal}. The apparent sign difference comes solely from the convention $D_j=\partial_j$ fixed at the beginning of the chapter.

\begin{corollary}[The ideal of smoothing operators]
\label{cor:ideal-regularizantes-pseudo}
Let $U\subseteq\mathbb R^n$ be open. If $R\in\Psi^{-\infty}(U;\mathbb C^r,\mathbb C^s)$, then $R_h^*\in\Psi^{-\infty}(U;\mathbb C^s,\mathbb C^r)$. If $A\colon C_c^\infty(U,\mathbb C^s)\to C^\infty(U,\mathbb C^q)$ and $B\colon C_c^\infty(U,\mathbb C^p)\to C^\infty(U,\mathbb C^r)$ are pseudodifferential operators with support conditions making the compositions well defined, then
\[
AR\in\Psi^{-\infty}_{\mathrm{loc}}(U;\mathbb C^r,\mathbb C^q),
\qquad
RB\in\Psi^{-\infty}_{\mathrm{loc}}(U;\mathbb C^p,\mathbb C^s).
\]
The compositions belong to the classes without the subscript $\mathrm{loc}$ when their kernels are properly supported; this holds, in particular, when the factors are properly supported.
\end{corollary}

\begin{proof}
At the symbol level, each term in the composition and adjoint expansions contains derivatives of a symbol of arbitrarily negative order. It therefore remains in $S^{-\infty}(U\times\mathbb R^n;
\operatorname{Hom}(\mathbb C^p,\mathbb C^q))$, with target spaces adjusted in each composition. The kernel formulation gives the same conclusion: composition integrates a smooth kernel against a kernel singular only on the diagonal, and proper support allows the integral to be taken and differentiated.
\end{proof}

\subsection{Coordinate changes}

The preceding formulas were obtained in Euclidean coordinates. For them to define a theory on manifolds, we must check that a change of chart transforms a pseudodifferential operator into another of the same order and that the principal component follows the expected cotangent rule. Let $\kappa\colon U\longrightarrow V$ be a diffeomorphism between Euclidean open sets, and let $T_\kappa u:=u\circ\kappa$. If $A=\operatorname{Op}(a)$ acts on $V$, consider
\[
\widetilde A:=T_\kappa A T_\kappa^{-1}.
\]
For $x,y\in U$ sufficiently close, define
\[
B(x,y):=\int_0^1D\kappa\bigl(y+t(x-y)\bigr)\,dt.
\]
Hadamard's lemma~\ref{lem:hadamard-finito-dimensional}, with base point $y$, gives the exact identity $\kappa(x)-\kappa(y)=B(x,y)(x-y)$ and also $B(x,x)=D\kappa(x)$. After the change of variables $z=\kappa(y)$ and $\xi=B(x,y)^T\eta$, the transformed operator has amplitude
\[
\widetilde b(x,y,\xi)
=
a\bigl(\kappa(x),B(x,y)^{-T}\xi\bigr)
\frac{|\det D\kappa(y)|}{|\det B(x,y)|}.
\]
Away from the diagonal the kernel is smooth, so we may insert a cutoff and assume that $B(x,y)$ is invertible.

\begin{theorem}[Invariance under coordinate changes]
\label{teo:cambio-coordenadas-pseudo}
\index{pseudodifferential operator!coordinate change}
The operator $\widetilde A$ is pseudodifferential of the same order as $A$. If $a$ is classical of order $m$, the transformed principal symbol is
\[
\widetilde a_m(x,\xi)
=
a_m\bigl(\kappa(x),D\kappa(x)^{-T}\xi\bigr).
\]
The full symbol is obtained by reducing $\widetilde b$ using Theorem~\ref{teo:amplitud-a-simbolo} and depends on higher derivatives of $\kappa$.
\end{theorem}

\begin{proof}
The preceding formula shows that $\widetilde b$ is an amplitude of order $m$ on each pair of compact coordinate sets. This follows from the chain rule: each derivative in $\xi$ produces a factor $B^{-T}$ and a derivative in the covariable of $a$, hence lowers the order; derivatives in $x$ and $y$ produce only smooth coefficients and spatial derivatives of $a$. Theorem~\ref{teo:amplitud-a-simbolo} produces the full symbol; weights in that reduction are controlled by Peetre's inequality in Lemma~\ref{lem:desigualdad-peetre-peso-japones}. At $y=x$, the ratio of determinants equals one and $B(x,x)=D\kappa(x)$. The $\alpha=0$ term of the reduction gives the principal formula, while each term with $|\alpha|\geq1$ contains a derivative in $\xi$ and has order at most $m-1$.
\end{proof}

This formula proves the two assertions needed globally. The full symbol is not an intrinsic function, because it changes by terms involving higher derivatives of the diffeomorphism. The principal component, however, transforms exactly by the cotangent action: if $\eta\in T_{\kappa(x)}^*V$ and $\xi=D\kappa(x)^T\eta$, then $\eta=D\kappa(x)^{-T}\xi$. It therefore defines a global object on the cotangent bundle.

Amplitude reduction controls the composition and adjoint remainders, while coordinate changes distinguish the chart-dependent full symbol from its cotangent principal component. This organization is presented in detail in \cite[Chapter~II, Sections~3--5]{Taylor1974Pseudo} and, in a later, more compact treatment, in \cite[Sections~0.2--0.4]{TaylorNonlinearPseudo}.

\section{Pseudodifferential operators on manifolds}
\label{sec:pseudo-variedades}

Invariance under coordinate changes allows us to leave the Euclidean model. The global definition will be local near the diagonal, supplemented by a proper-support condition when we wish to act on distributions without compact support.

Throughout this section, $M$ is a finite-dimensional smooth manifold without boundary. It is not assumed compact, complete, or Riemannian until the subsection on geometric quantization. Let $\mathbf{E},\mathbf{F}\to M$ be smooth complex vector bundles of finite rank. The real case follows by complexification and restriction.

\subsection{Local definition and independence of choices}

Let $(U,\kappa)$ be a chart, and let $\mathbf{e}=(\mathbf{e}_1,\dots,\mathbf{e}_r)$ and $\mathbf{f}=(\mathbf{f}_1,\dots,\mathbf{f}_s)$ be local frames of $\mathbf{E}$ and $\mathbf{F}$. If $\chi,\psi\in C_c^\infty(U)$, the localized operator $\chi A\psi$ is represented, using the chart and frames, by a matrix of operators between functions compactly supported in $\kappa(U)$. Extension by zero after shrinking the supports allows it to be regarded on $\mathbb R^n$.

\begin{definition}[Pseudodifferential operator on a manifold]
\label{def:pseudo-variedad-local}
\index{pseudodifferential operator!on a manifold}
\glsadd{operador-pseudodiferencial}
Let $M$ be a finite-dimensional smooth manifold without boundary, and let $\mathbf{E},\mathbf{F}\to M$ be smooth complex vector bundles of finite rank. For $m\in\mathbb R$, a continuous linear operator \[
A\colon \Gamma_c(M,\mathbf{E})\longrightarrow \Gamma(M,\mathbf{F})
\] belongs to $\Psi^m_{\mathrm{loc}}(M;\mathbf{E},\mathbf{F})$ if the following two conditions hold:
\begin{enumerate}[label=(\alph*)]
\item the kernel $\mathbf{K}_A$ is smooth off $\Delta_M$;
\item for every chart $(U,\kappa)$, every pair of local frames of $\mathbf{E}$ and $\mathbf{F}$ over $U$, and all $\chi,\psi\in C_c^\infty(U)$, each entry of the Euclidean representation of $\chi A\psi$ is $\operatorname{Op}(a_{ba})$ with $a_{ba}\in S^m_{1,0}(\kappa(U))$, modulo an operator with smooth kernel.
\end{enumerate}
If the representations can be chosen classical, write $\Psi^m_{\mathrm{cl,loc}}(M;\mathbf{E},\mathbf{F})$. Denote by \[
\Psi^m(M;\mathbf{E},\mathbf{F})\subseteq\Psi^m_{\mathrm{loc}}(M;\mathbf{E},\mathbf{F})
\] the subclass of properly supported operators. This latter notation will be used for composition and global action on distributions.
\end{definition}

\begin{definition}[Smoothing operators on a manifold]
\label{def:operadores-regularizantes-variedad}
\index{smoothing operator!on a manifold}
Let $M$ be a finite-dimensional smooth manifold without boundary, and let $\mathbf{E},\mathbf{F}\to M$ be smooth complex vector bundles of finite rank. Define
\[
\Psi^{-\infty}_{\mathrm{loc}}(M;\mathbf{E},\mathbf{F})
:=
\left\{
A\in\Psi^m_{\mathrm{loc}}(M;\mathbf{E},\mathbf{F})\text{ for some }m\in\mathbb R
\;\middle|\;
\mathbf{K}_A\text{ is smooth on }M\times M
\right\}.
\]
Proposition~\ref{prop:equiv-simbolo-nucleo-regularizante}, applied in charts and frames, shows that an element of this class locally belongs to every order. The properly supported class is
\begin{equation}
\Psi^{-\infty}(M;\mathbf{E},\mathbf{F})
:=
\left\{
A\in\Psi^{-\infty}_{\mathrm{loc}}(M;\mathbf{E},\mathbf{F})
\;\middle|\;
\mathbf{K}_A\text{ is properly supported}
\right\}.
\label{eq:def-operadores-regularizantes-variedad}
\end{equation}
If $M$ is closed, the two classes coincide. If $M$ is noncompact, smoothness of the kernel implies neither proper support nor, after a Riemannian metric and Hermitian bundle metrics are chosen, compactness of the operator $L^2(M,\mathbf{E})\longrightarrow L^2(M,\mathbf{F})$.
\end{definition}

The definition uses charts, frames, and cutoffs, but these choices merely express the kernel locally. The following result shows that it suffices to verify the definition in a trivializing atlas and that the resulting class is intrinsic.

\begin{lemma}[A properly supported cutoff near the diagonal]
\label{lem:corte-propio-diagonal-pseudo}
Let $M$ be a finite-dimensional smooth manifold without boundary. There is
$\rho\in C^\infty(M\times M,[0,1])$ equal to one on a neighborhood of
$\Delta_M$, whose support is proper with respect to both projections to $M$.
\end{lemma}

\begin{proof}
Choose a countable locally finite cover $(U_j)_{j\in\mathbb N}$ by
relatively compact coordinate domains, and a smooth partition of unity
$(\eta_j)_{j\in\mathbb N}$ with
$\operatorname{supp}\eta_j\Subset U_j$. Set
\[
 b(x):=\left(\sum_{j=1}^{\infty}\eta_j(x)^2\right)^{\frac{1}{2}},
 \qquad \theta_j(x):=\frac{\eta_j(x)}{b(x)}.
\]
The sum is locally finite and $b(x)>0$, so $b$ and the functions $\theta_j$
are smooth. Moreover, $\sum_j\theta_j(x)^2=1$. Define
\[
 r(x,y):=\sum_{j=1}^{\infty}\theta_j(x)\theta_j(y).
\]
This sum is also locally finite; Cauchy--Schwarz gives
$0\leq r(x,y)\leq1$, and $r(x,x)=1$. Take
$\chi\in C^\infty(\mathbb R,[0,1])$ equal to zero on
$(-\infty,\displaystyle\frac{1}{4}]$ and equal to one on $[\displaystyle\frac{1}{2},+\infty)$. Then
$\rho:=\chi\circ r$ equals one on the open set $\{r>\displaystyle\frac{1}{2}\}$, which
contains the diagonal.

Write $K_j:=\operatorname{supp}\theta_j$. The family $(K_j)$ is locally
finite, and every $K_j$ is compact. Thus
$S:=\bigcup_j(K_j\times K_j)$ is closed in $M\times M$.
Outside $S$ we have $r=0$, so $\operatorname{supp}\rho\subseteq S$.
For each compact $K\subseteq M$, only finitely many $K_j$ meet $K$.
Consequently, $\operatorname{supp}\rho\cap(K\times M)$ is a closed subset
of the finite union of the compact sets $(K\cap K_j)\times K_j$, and is
therefore compact. Interchanging the variables proves that
$\operatorname{supp}\rho\cap(M\times K)$ is compact as well.
These are precisely the two proper support conditions.
\end{proof}

\begin{theorem}[Independence of charts, frames, and cutoffs]
\label{teo:independencia-definicion-pseudo-variedad}
Let $M$ be a finite-dimensional smooth manifold without boundary, let $\mathbf{E},\mathbf{F}\to M$ be smooth complex vector bundles of finite rank, and let $m\in\mathbb R$. Definition~\ref{def:pseudo-variedad-local} is independent of the atlas, local frames, and cutoffs. Moreover, every $A\in\Psi^m_{\mathrm{loc}}(M;\mathbf{E},\mathbf{F})$ can be written as $A=A_0+R$, where $A_0\in\Psi^m(M;\mathbf{E},\mathbf{F})$ and $R\in\Psi^{-\infty}_{\mathrm{loc}}(M;\mathbf{E},\mathbf{F})$.
\end{theorem}

\begin{proof}
Consider two overlapping charts. After inserting cutoffs whose supports remain inside the overlap, the representations are conjugate by the transition diffeomorphism. Theorem~\ref{teo:cambio-coordenadas-pseudo} shows that conjugation preserves the order and the classical class. A change of frame in $\mathbf{E}$ or $\mathbf{F}$ multiplies the operator on the right or left by a smooth matrix. Theorem~\ref{teo:composicion-simbolica-pseudo}, applied to symbols independent of $\xi$, again shows that the order is preserved. Changing a cutoff yields a difference supported where the two variables are distinct, or another smooth multiplication; in the first case the kernel is smooth, and in the second the composition calculus applies.

For the last assertion, take the cutoff $\rho$ from
Lemma~\ref{lem:corte-propio-diagonal-pseudo}. The kernel
$\rho\mathbf K_A$ defines a properly supported operator $A_0$.
Since $\rho=1$ near the diagonal, the difference
$(1-\rho)\mathbf K_A$ is smooth on all of $M\times M$. Thus
$R:=A-A_0\in\Psi^{-\infty}_{\mathrm{loc}}(M;\mathbf E,\mathbf F)$.
In each chart, reduction of the amplitude obtained by multiplication
by $\rho$ preserves order $m$, so $A_0\in\Psi^m(M;\mathbf E,\mathbf F)$.
\end{proof}

The distinction between the two classes is not merely typographical. On a noncompact manifold, the purely local definition does not control kernel behavior when both variables leave every compact set. Without proper support, a composition may be undefined even when both factors are locally pseudodifferential. On a closed manifold, every projection of a closed subset of $M\times M$ is proper, so the distinction disappears.

\subsection{The global principal symbol}

Once independence of local choices has been established, the principal homogeneous components can be glued over the cotangent bundle. The result is the global object that will preserve composition, adjoints, and ellipticity.

Let $\pi\colon T^*M\longrightarrow M$ be the projection and $0_M$ the zero section. In a chart and local frames, a classical operator $A\in\Psi^m_{\mathrm{cl,loc}}(M;\mathbf{E},\mathbf{F})$ has a principal component $a_m(x,\xi)$. Theorem~\ref{teo:cambio-coordenadas-pseudo} and the frame transformation show that, on overlaps, these components satisfy exactly the transformation rule for a bundle homomorphism \[
\pi^*\mathbf{E}\longrightarrow\pi^*\mathbf{F}
\] over $T^*M\setminus0_M$.

\begin{definition}[Pseudodifferential principal symbol]
\label{def:simbolo-principal-pseudodiferencial-global}
\index{principal symbol@principal symbol!pseudodifferential}
Let $M$ be a finite-dimensional smooth manifold without boundary, let $\mathbf{E},\mathbf{F}\to M$ be smooth complex vector bundles of finite rank, and let $\pi\colon T^*M\to M$. For $A\in\Psi^m_{\mathrm{cl,loc}}(M;\mathbf{E},\mathbf{F})$, the \textbf{principal Fourier symbol} of $A$ is the smooth section
\[
\boldsymbol{\sigma}_m^\Psi(A)
\in
\Gamma\bigl(T^*M\setminus0_M,
\operatorname{Hom}(\pi^*\mathbf{E},\pi^*\mathbf{F})\bigr)
\]
whose local expression is $a_m$. It satisfies
\[
\boldsymbol{\sigma}_m^\Psi(A)(x,t\xi)=t^m\boldsymbol{\sigma}_m^\Psi(A)(x,\xi),
\qquad t>0.
\]
\end{definition}

For a differential operator $P$ of order $m$, the exact compatibility with Chapter~\ref{cap:simbolo-principal-operadores-elipticos} is
\begin{equation}
\label{eq:conversion-simbolos-diferencial-pseudo}
\boldsymbol{\sigma}_m^\Psi(P)(x,\xi)
=
i^m\boldsymbol{\sigma}_m(P)(x,\xi).
\end{equation}
The factor $i^m$ comes from the convention $D_j=\partial_j$. Since this factor is a nonzero scalar, the two notions of ellipticity coincide. The earlier composition and adjoint rules are also recovered exactly, as checked after Theorem~\ref{teo:simbolo-adjunto-pseudo}.

\begin{example}[The Dirac operator as a pseudodifferential operator]
\label{ej:dirac-como-pseudodiferencial}
\index{Dirac operator!as a pseudodifferential operator}
Let \(M\) be a Riemannian manifold without boundary, and let \(D=c\circ\nabla^{\mathbf{S}}\) be a Dirac-type operator constructed in Section~\ref{sec:operadores-dirac-diferenciales}. Every differential operator of order one is a properly supported classical pseudodifferential operator; thus
\[
 D\in\Psi^1_{\mathrm{cl}}(M;\mathbf{S},\mathbf{S}).
\]
Combining \eqref{eq:conversion-simbolos-diferencial-pseudo} with \eqref{eq:simbolo-diferencial-dirac} gives, with this chapter's Fourier convention,
\begin{equation}
 \boldsymbol\sigma_1^\Psi(D)(x,\xi)=i\,c(\xi).
\label{eq:simbolo-pseudodiferencial-dirac}
\end{equation}
For \(\xi\neq0\), this homomorphism is invertible and
\[
 \bigl(\boldsymbol\sigma_1^\Psi(D)(x,\xi)\bigr)^{-1}
 =\frac{i\,c(\xi)}{|\xi|_{\mathbf{g}}^2},
\]
since \((i\,c(\xi))^2=|\xi|_{\mathbf{g}}^2I\). If \(\mathbf{S}\) is graded, the same formula restricted to \(\mathbf{S}_x^+\to\mathbf{S}_x^-\) is the pseudodifferential symbol of \(D^+\).
\end{example}

\begin{definition}[Space of principal symbols]
\label{def:espacio-simbolos-principales-globales}
Let $M$ be a finite-dimensional smooth manifold without boundary, let $\mathbf{E},\mathbf{F}\to M$ be smooth complex vector bundles of finite rank, and let $\pi\colon T^*M\to M$. Let $\mathcal S^m_{\mathrm{hom}}(T^*M;\operatorname{Hom}(\mathbf{E},\mathbf{F}))$ be the space of smooth sections of $\operatorname{Hom}(\pi^*\mathbf{E},\pi^*\mathbf{F})$ over $T^*M\setminus0_M$ that are homogeneous of degree $m$ under positive fiber dilations.
\end{definition}

\begin{theorem}[Exact sequence of the principal symbol]
\label{teo:sucesion-exacta-simbolo-principal-pseudo}
Let $M$ be a finite-dimensional smooth manifold without boundary, and let $\mathbf{E},\mathbf{F},\mathbf{G}\to M$ be smooth complex vector bundles of finite rank. For every $m\in\mathbb R$, there is an exact sequence
\[
0\longrightarrow
\Psi^{m-1}_{\mathrm{cl}}(M;\mathbf{E},\mathbf{F})
\longrightarrow
\Psi^m_{\mathrm{cl}}(M;\mathbf{E},\mathbf{F})
\xrightarrow{\ \boldsymbol{\sigma}_m^\Psi\ }
\mathcal S^m_{\mathrm{hom}}(T^*M;\operatorname{Hom}(\mathbf{E},\mathbf{F}))
\longrightarrow0.
\]
Moreover, for every $\ell\in\mathbb R$, if $A\in\Psi^m_{\mathrm{cl}}(M;\mathbf{E},\mathbf{F})$ and $B\in\Psi^\ell_{\mathrm{cl}}(M;\mathbf{G},\mathbf{E})$, then
\[
\boldsymbol{\sigma}_{m+\ell}^\Psi(AB)
=
\boldsymbol{\sigma}_m^\Psi(A)\circ\boldsymbol{\sigma}_\ell^\Psi(B)
\]
and, when Hermitian bundle metrics and a density have been fixed to define the adjoint,
\[
\boldsymbol{\sigma}_m^\Psi(A_h^*)=\boldsymbol{\sigma}_m^\Psi(A)^*.
\]
\end{theorem}

\begin{proof}
Well-definedness and the two compatibility statements follow from Theorems~\ref{teo:cambio-coordenadas-pseudo}, \ref{teo:composicion-simbolica-pseudo}, and~\ref{teo:simbolo-adjunto-pseudo}. If $A$ has order at most $m-1$, its homogeneous component of degree $m$ vanishes. Conversely, if this component vanishes in every chart, the definition of a classical symbol shows that each local representative satisfies estimates of order $m-1$; hence $A\in\Psi^{m-1}_{\mathrm{cl}}(M;\mathbf{E},\mathbf{F})$.

For surjectivity, let $s$ be a section of the homogeneous symbol space
in the statement. Choose a countable locally finite cover $(U_j)$ by
coordinate domains trivializing both bundles, and smooth functions
$\theta_j$ compactly supported in $U_j$ with $\sum_j\theta_j^2=1$, as in
Lemma~\ref{lem:corte-propio-diagonal-pseudo}. Let $s_j(x,\xi)$ be the
matrix representing $s$ in the corresponding coordinates and frames.
Choose $\zeta_j\in C_c^\infty(U_j)$ equal to one near
$\operatorname{supp}\theta_j$, and a smooth cutoff $c(\xi)$ equal to zero
for $|\xi|\leq1$ and equal to one for $|\xi|\geq2$.
The matrix $a_j(x,\xi):=\zeta_j(x)c(\xi)s_j(x,\xi)$, extended by zero
in $x$ outside the chart and in $\xi$ near zero, is a classical symbol
of order $m$. Indeed, homogeneity gives
\[
 |\partial_x^\alpha\partial_\xi^\beta s_j(x,\xi)|
 \leq C_{\alpha,\beta,j}|\xi|^{m-|\beta|}
 \quad(x\in\operatorname{supp}\zeta_j,\ |\xi|\geq1),
\]
and derivatives of the cutoff in $\xi$ only produce terms compactly
supported in that variable.

Let $A_j$ be the local quantization of $a_j$, transported to the bundles.
Define
\[
 A\mathbf u:=\sum_{j=1}^{\infty}\theta_j A_j(\theta_j\mathbf u),
 \qquad \mathbf u\in\Gamma_c(M,\mathbf E).
\]
Only finitely many indices contribute for a compactly supported section.
The kernel of each summand is supported in
$\operatorname{supp}\theta_j\times\operatorname{supp}\theta_j$.
Local finiteness and compactness of these supports show, as in the cited
lemma, that the kernel of $A$ is proper in both variables.
Each compact localization involves a finite sum, so the composition and
coordinate change calculus gives
$A\in\Psi^m_{\mathrm{cl}}(M;\mathbf E,\mathbf F)$.
The principal component of the $j$th summand is $\theta_j^2s$, since
$\zeta_j=1$ near its support. Hence
\[
 \boldsymbol\sigma_m^\Psi(A)=\sum_j\theta_j^2s=s.
\]
This proves surjectivity.
\end{proof}

\subsection{The conormal bundle of the diagonal}

The definition by symbols still depends on a local quantization. There is an equivalent description using only the kernel: classical pseudodifferential operators are exactly those whose kernels are conormal to the diagonal. To formulate this equivalence, we first fix the order and density conventions for a general submanifold.

Let $X$ be a smooth manifold without boundary of dimension $d$, let $Y\subseteq X$ be an embedded submanifold without boundary of codimension $k$, and let $\mathbf{G}\to X$ be a smooth complex vector bundle of finite rank $r_{\mathbf{G}}$. Its normal and conormal bundles are
\[
NY:=\operatorname{coker}\bigl(TY\hookrightarrow TX|_Y\bigr),
\qquad
N^*Y:=\{(y,\eta)\in T^*X|_Y\mid \eta|_{T_yY}=0\}.
\]
Let $\pi_Y\colon N^*Y\to Y$ denote the projection and $\boldsymbol{\mathcal D}_1(Y)$ the bundle of weight-one densities from Definition~\ref{def:variedades-riemannianas-haz-de-densidades-de-peso-w-en-una-variedad}. A chart $(y,z)\in\mathbb R^{d-k}\times\mathbb R^k$ is \textbf{adapted to $Y$} if $Y$ is given by $z=0$. In this chart, conormal covectors are precisely $(y,0;0,\theta)$.

\begin{definition}[Conormal distribution]
\label{def:distribucion-conormal-general}
\index{conormal distribution@conormal distribution}
A distribution $u\in\mathcal D'(X,\mathbf{G})$ is \textbf{classical conormal to $Y$ of order $m$}, written
\[
u\in I^m(X,Y;\mathbf{G}),
\]
if it is smooth on $X\setminus Y$ and, in every localization in an adapted chart and a trivialization of $G$, can be written modulo a smooth section as
\begin{equation}
u
=
\left[
(2\pi)^{-k}\operatorname{Os}\!\int_{\mathbb R_\theta^k}
e^{iz\cdot\theta}b(y,z,\theta)\,d\theta
\right]|dy\,dz|,
\qquad
b\in
S^{m+\frac{d-2k}{4}}_{\mathrm{cl}}
\bigl(\mathbb R^d_{(y,z)}\times\mathbb R^k_\theta;\mathbb C^{r_{\mathbf{G}}}\bigr).
\label{eq:representacion-conormal-general}
\end{equation}
Here $|dy\,dz|$ records the coordinate density implicit in the convention defining distributions as functionals on test sections of the dual bundle. This factor is suppressed in Euclidean formulas but must be retained when analyzing coordinate changes. After localization, the amplitude is taken compactly supported in the spatial variables. The number $m$ is the \textbf{conormal order}; the shift $\displaystyle\frac{d-2k}{4}$ adapts Hörmander's convention to the operator order used on the diagonal.
\end{definition}

The definition is independent of adapted coordinates, as the following calculation shows. Let \[
\widetilde y=\phi(y,z),
\qquad
\widetilde z=\psi(y,z)
\] be two adapted systems related by the diffeomorphism $\Phi=(\phi,\psi)$. The identity $\psi(y,0)=0$ and Hadamard's lemma~\ref{lem:hadamard-finito-dimensional}, applied with $y$ as parameter, give the exact factorization
\begin{equation}
\psi_a(y,z)=\sum_{b=1}^k C_{ab}(y,z)z_b,
\qquad
C_{ab}(y,z):=\int_0^1
\partial_{z_b}\psi_a(y,tz)\,dt.
\label{eq:cambio-normal-hadamard-conormal}
\end{equation}
In particular, $C(y,0)=D_z\psi(y,0)$. Since \[
D\Phi(y,0)=
\begin{pmatrix}
D_y\phi(y,0)&D_z\phi(y,0)\\
0&C(y,0)
\end{pmatrix}
\] is invertible, so are $A(y):=D_y\phi(y,0)$ and $C(y,0)$. After shrinking the chart, $C(y,z)$ remains invertible. For each fixed $(y,z)$, the linear substitution
\begin{equation}
\theta=C(y,z)^T\widetilde\theta,
\qquad
d\theta=|\det C(y,z)|\,d\widetilde\theta
\label{eq:cambio-covariable-conormal}
\end{equation}
is exact and satisfies $z\cdot\theta=\widetilde z\cdot\widetilde\theta$. If $T(y,z)$ denotes the matrix taking the components of $G$ in the original frame to those in the new frame, the transformed amplitude is
\begin{equation}
\begin{split}
\widetilde b(\widetilde y,\widetilde z,\widetilde\theta)
:={}&T(y,z)b\bigl(y,z,C(y,z)^T\widetilde\theta\bigr)
|\det C(y,z)||\det D\Phi(y,z)|^{-1},
\qquad (y,z)=\Phi^{-1}(\widetilde y,\widetilde z).
\end{split}
\label{eq:amplitud-transformada-conormal}
\end{equation}
Indeed, the last determinant is exactly the one that converts $|dy\,dz|$ into $|d\widetilde y\,d\widetilde z|$. To check the order, first apply derivatives in $\widetilde\theta$ and then spatial derivatives. The Leibniz and chain rules express each derivative of $b\bigl(\Phi^{-1}(\widetilde y,\widetilde z),
C^T\widetilde\theta\bigr)$ as a finite sum of terms
\[
A_{\alpha,\beta,\gamma;\delta,\nu,\rho}
(\widetilde y,\widetilde z)\widetilde\theta^\rho
(D_{y,z}^\delta D_\theta^\nu b)
\bigl(y,z,C(y,z)^T\widetilde\theta\bigr),
\]
where $(y,z)=\Phi^{-1}(\widetilde y,\widetilde z)$ and
\[
|\nu|=|\gamma|+|\rho|,
\qquad
|\delta|+|\rho|\leq|\alpha|+|\beta|.
\]
Indeed, a spatial derivative hitting $b$ increases $|\delta|$; if it hits $C^T\widetilde\theta$, it creates both a linear factor in $\widetilde\theta$ and a derivative in $\theta$; if it hits a coefficient, it merely differentiates that coefficient. This proves both relations by induction on $|\alpha|+|\beta|$. Uniform invertibility of $C$ on each localized compact set gives $c\langle\widetilde\theta\rangle\leq
\langle C^T\widetilde\theta\rangle\leq
C\langle\widetilde\theta\rangle$. The factor $\widetilde\theta^\rho$ is offset by the $|\rho|$ additional derivatives in $\theta$. Since the remaining factors in \eqref{eq:amplitud-transformada-conormal} are smooth, for all multi-indices $\alpha,\beta,\gamma$ there is a constant $C_{\alpha,\beta,\gamma}>0$ such that
\[
\bigl\|D_{\widetilde y}^{\alpha}D_{\widetilde z}^{\beta}
D_{\widetilde\theta}^{\gamma}\widetilde b\bigr\|
\leq C_{\alpha,\beta,\gamma}
\langle\widetilde\theta\rangle^{m+\frac{d-2k}{4}-|\gamma|}.
\]
Each homogeneous component transforms by the same formula; thus \eqref{eq:amplitud-transformada-conormal} belongs to the same classical class. Applying the calculation to $\Phi^{-1}$ proves the reverse implication.

The principal coefficient alone is not a canonical function. If $\mu=m+\displaystyle\frac{d-2k}{4}$ and $b_\mu$ is the principal homogeneous component, then
\begin{equation}
b_\mu(y,0,\theta)|dy|
\label{eq:simbolo-conormal-con-densidad}
\end{equation}
defines a homogeneous section of degree $\mu$ of
\[
\pi_Y^*\bigl(G\restriction_Y\otimes\boldsymbol{\mathcal D}_1(Y)\bigr)
\quad\text{on }N^*Y\setminus0_Y,
\]
modulo the next degree. This is the \textbf{principal conormal symbol} in the density convention fixed above. To verify invariance, restrict \eqref{eq:amplitud-transformada-conormal} to $z=0$. The triangular form of $D\Phi(y,0)$ implies \[
|\det D\Phi(y,0)|=|\det A(y)|\,|\det C(y,0)|,
\] and hence
\begin{equation}
\widetilde b_\mu
\bigl(\phi(y,0),0,\widetilde\theta\bigr)
=T(y,0)b_\mu\bigl(y,0,C(y,0)^T\widetilde\theta\bigr)
|\det A(y)|^{-1}
\quad\bmod\{\text{conormal symbols of degree at most }\mu-1\}.
\label{eq:transformacion-simbolo-conormal-principal}
\end{equation}
Since $|d\widetilde y|=|\det A(y)|\,|dy|$ on $Y$, multiplying \eqref{eq:transformacion-simbolo-conormal-principal} by $|d\widetilde y|$ cancels the sole remaining determinant. In turn, $\theta=C(y,0)^T\widetilde\theta$ is exactly the coordinate transformation rule in $N^*Y$. Thus $b_\mu(y,0,\theta)|dy|$ defines the asserted global section. This verification explains why omitting the density line would produce a noncanonical object.

The normal phases used in adapted charts are related by the homogeneous linear substitution \eqref{eq:cambio-covariable-conormal}. Since no additional phase variables are introduced, the changes are fully described by the cotangent rule and the preceding densities. With this convention, the normalization agrees with \cite[Theorems~18.2.8--18.2.11]{HormanderIII}: in ambient dimension $d$ and codimension $k$, a conormal distribution of order $m$ is represented by an amplitude of order $m+\displaystyle\frac{d-2k}{4}$.

For $X=M\times M$ and $Y=\Delta_M$, the tangent space at $(x,x)$ consists of pairs $(v,v)$. A covector $(\xi,\eta)$ annihilates it if and only if $\xi(v)+\eta(v)=0$ for every $v\in T_xM$. Thus
\begin{equation}
\label{eq:conormal-diagonal-identificacion}
N^*\Delta_M
=
\{(x,\xi;x,-\xi)\mid x\in M,\ \xi\in T_x^*M\}.
\end{equation}
This is the geometric reason that the cotangent bundle, rather than the tangent bundle, parametrizes the kernel's oscillations.

In this case $d=2n$ and $k=n$, so the shift in Definition~\ref{def:distribucion-conormal-general} is $\displaystyle\frac{2n-2n}{4}=0$. A distribution $\mathbf{K}\in I^m(M\times M,\Delta_M;\mathbf{F}\boxtimes \mathbf{E}^*)$ therefore has the local expression, modulo a smooth section,
\[
(2\pi)^{-n}\operatorname{Os}\!\int_{\mathbb R_\xi^n}
e^{i(x-y)\cdot\xi}b(x,y,\xi)\,d\xi,
\qquad b\text{ classical amplitude of order }m.
\]
The bundle $\boldsymbol{\mathcal D}_1$ of Definition~\ref{def:variedades-riemannianas-haz-de-densidades-de-peso-w-en-una-variedad} supplies the required density factor. To relate it to an operator kernel, temporarily fix a positive density $\nu$ on $M$ and use $\nu$ both to identify smooth sections with distributions and to write integration in the input variable. The principal conormal symbol of the kernel on the diagonal then takes values in
\[
\operatorname{Hom}(\mathbf{E},\mathbf{F})\otimes\boldsymbol{\mathcal D}_1(\Delta_M).
\]
The identification $\Delta_M\cong M$ and division by $\nu$ produce a section of $\operatorname{Hom}(\pi^*\mathbf{E},\pi^*\mathbf{F})$. If $\widetilde\nu=f\nu$, the kernel coefficient relative to integration in the input variable acquires $f(y)^{-1}$, which cancels with the new integration density $f(y)\nu(y)$. Pairing in the output variable, in contrast, multiplies the kernel distribution by $f(x)$; its conormal symbol on the diagonal therefore acquires $f\restriction_{\Delta_M}$, and the final division by $\widetilde\nu=f\nu$ cancels exactly this factor. The normalized section is consequently independent of $\nu$. In a chart, taking the coordinate density, this section is precisely the homogeneous component of the Kohn--Nirenberg amplitude. Thus normalization by the density preserves the operator order and produces an intrinsic principal conormal symbol.

\begin{theorem}[Characterization by conormal kernels]
\label{teo:equiv-pseudo-nucleo-conormal}
Let $M$ be a finite-dimensional smooth manifold without boundary, and let $\mathbf{E},\mathbf{F}\to M$ be smooth complex vector bundles of finite rank. Let \[
A\colon\Gamma_c(M,\mathbf{E})\longrightarrow\Gamma(M,\mathbf{F})
\] be a continuous properly supported operator. The following conditions are equivalent:
\begin{enumerate}[label=(\alph*)]
\item $A\in\Psi^m_{\mathrm{cl}}(M;\mathbf{E},\mathbf{F})$;
\item $\mathbf{K}_A\in I^m(M\times M,\Delta_M;\mathbf{F}\boxtimes \mathbf{E}^*)$.
\end{enumerate}
The equivalence preserves the principal symbol: the homogeneous component of the conormal amplitude, restricted to $N^*\Delta_M$ by \eqref{eq:conormal-diagonal-identificacion} and normalized by the density as explained above, is $\boldsymbol{\sigma}_m^\Psi(A)$.
\end{theorem}

\begin{proof}
Suppose (a) holds. In a chart and local frames, the kernel formula for $\operatorname{Op}(a)$ is an oscillatory representation with amplitude $a(x,\xi)$. Chart changes transform the phase and amplitude as in the proof of Theorem~\ref{teo:cambio-coordenadas-pseudo}; the principal component follows the cotangent rule. Off the diagonal, the integrations by parts in the first subsection of \ref{sec:pseudo-nucleos-calculo} yield smoothness. Thus (b) holds.

Conversely, suppose (b) holds. Each local conormal representation defines an operator with amplitude of order $m$. Theorem~\ref{teo:amplitud-a-simbolo} reduces it to a classical symbol of that order. The remaining smooth term is smoothing by Proposition~\ref{prop:equiv-simbolo-nucleo-regularizante}. This gives (a). In both directions, the $\alpha=0$ term of amplitude reduction preserves the principal homogeneous component; the other terms contain at least one derivative in $\xi$ and lower the order. This proves the last assertion.
\end{proof}

Combining the preceding theorem with Theorem~\ref{teo:amplitud-a-simbolo} gives a precise equivalence among the three descriptions:
\[
\begin{aligned}
\text{classical symbol of order }m
&\quad\Longleftrightarrow\quad
\text{classical amplitude of order }m,\\
\text{classical amplitude of order }m
&\quad\Longleftrightarrow\quad
\text{conormal kernel of order }m\text{ to }\Delta_M.
\end{aligned}
\]
The first equivalence is exact at the operator level and determines the full symbols, in each chart and frames, modulo the corresponding class $S^{-\infty}(U\times\mathbb R^n;
\operatorname{Hom}(\mathbb C^r,\mathbb C^s))$; the second is exact modulo smooth kernels. The principal symbol is exactly the same global object in all three formulations.

\section{Riemannian quantization via the exponential map}
\label{sec:pseudo-cuantizacion-geometrica}

The preceding global definition is obtained by localizing Kohn--Nirenberg quantization. Once a metric and bundle connections are fixed, one can also construct a geometric quantization: the difference $x-y$ is replaced by the vector $\exp_x^{-1}(y)$, and fibers are compared by parallel transport. This formulation does not remove coordinates from the proofs, but clearly separates the geometric choices from the symbolic part.

Now suppose that $(M,\mathbf{g})$ is Riemannian and without boundary, but not necessarily complete or compact. Let $\mathbf{E},\mathbf{F}$ be complex vector bundles with Hermitian bundle metrics and compatible connections. There is an open neighborhood $\mathcal V$ of the zero section of $TM$ on which
\[
(x,v)\longmapsto (x,\exp_xv)
\]
is a diffeomorphism onto a neighborhood $\mathcal U$ of $\Delta_M$. This assertion is local around each $(x,0)$ and is globalized by shrinking a locally finite cover; it requires no uniform injectivity radius. Write
\[
\ell(x,y):=\exp_x^{-1}(y)\in T_xM,
\qquad (x,y)\in\mathcal U.
\]

Let $\tau^{\mathbf{E}}_{x\leftarrow y}\colon \mathbf{E}_y\longrightarrow \mathbf{E}_x$ be parallel transport along the radial geodesic $t\mapsto\exp_x(t\ell(x,y))$. Choose a properly supported $\chi\in C^\infty(M\times M)$, supported in $\mathcal U$ and equal to one near the diagonal. The metric induces Lebesgue measures $dv_x$ on $T_xM$ and $d\xi_x$ on $T_x^*M$.

By the construction in Subsection~\ref{subsec:descomposicion-horizontal-vertical-conexion}, the Levi--Civita connection determines the horizontal--vertical decomposition of $T(T^*M)$, while the connections on $\mathbf{E}$ and $\mathbf{F}$ induce the pullback connection on $\pi^*\operatorname{Hom}(\mathbf{E},\mathbf{F})\to T^*M$. Differentiating with this connection in horizontal and vertical directions gives the derivatives $\nabla^{\mathrm H}$ and $\nabla^{\mathrm V}$ defined in \eqref{eq:definicion-derivadas-horizontal-vertical-cotangente}. Write $\langle\xi\rangle_{\mathbf{g}}=(1+\|\xi\|_{\mathbf{g}}^2)^{\displaystyle\frac{1}{2}}$.

Fix a chart $(U,x^1,\ldots,x^n)$, the induced coordinates $(x,\xi)$ determined by $\xi=\xi_jdx^j$, and frames $\mathbf{e}=(\mathbf{e}_1,\ldots,\mathbf{e}_r)$ and $\mathbf{f}=(\mathbf{f}_1,\ldots,\mathbf{f}_s)$ of $\mathbf{E}$ and $\mathbf{F}$. Adopt the conventions
\[
\nabla^{\mathbf{E}}_{\partial_i}\mathbf{e}_\alpha
=\sum_{\gamma=1}^{r} \mathbf{e}_\gamma(\omega_i^{\mathbf{E}})_{\gamma\alpha},
\qquad
\nabla^{\mathbf{F}}_{\partial_i}\mathbf{f}_b
=\sum_{c=1}^{s} \mathbf{f}_c(\omega_i^{\mathbf{F}})_{cb},
\qquad
\nabla_{\boldsymbol{\partial}_i}\boldsymbol{\partial}_j=\Gamma^k_{ij}\boldsymbol{\partial}_k.
\]
Proposition~\ref{prop:levantamiento-horizontal-cotangente} derived the sign in the formula
\[
H_i
=D_{x^i}+\sum_{j,k=1}^n\Gamma^k_{ij}(x)\xi_kD_{\xi_j};
\]
from the covector parallel-transport equation; this is identity \eqref{eq:levantamiento-horizontal-cotangente-local}. In particular, if $\mathbf{a}\in\Gamma(\pi^*\operatorname{Hom}(\mathbf{E},\mathbf{F})|_{T^*U})$ and we also denote its local matrix by $\mathbf{a}(x,\xi)$, Proposition~\ref{prop:formulas-locales-derivadas-horizontal-vertical} gives
\[
\nabla^{\mathrm V}_j \mathbf{a}=D_{\xi_j}\mathbf{a},
\qquad
\nabla^{\mathrm H}_i \mathbf{a}
=H_i \mathbf{a}+\omega_i^{\mathbf{F}} \mathbf{a}-\mathbf{a}\omega_i^{\mathbf{E}};
\]
these are identities \eqref{eq:derivadas-geometricas-cotangente-locales}.

To iterate these operations intrinsically, introduce the bundles
\[
\mathscr T^{p,q}
:=\pi^*\!\left(
(T^*M)^{\otimes p}\otimes(TM)^{\otimes q}
\otimes\operatorname{Hom}(\mathbf{E},\mathbf{F})
\right).
\]
The Levi--Civita connection on the tensor factors, the connection on $\operatorname{Hom}(\mathbf{E},\mathbf{F})$, and the horizontal--vertical decomposition induce
\[
\nabla^{\mathrm H}\colon\Gamma(\mathscr T^{p,q})
\longrightarrow\Gamma(\mathscr T^{p+1,q}),
\qquad
\nabla^{\mathrm V}\colon\Gamma(\mathscr T^{p,q})
\longrightarrow\Gamma(\mathscr T^{p,q+1}).
\]
As in Chapter~\ref{cap:sobolev-haces}, the new index is added at the beginning of the corresponding horizontal or vertical block. If $T$ has components $T_{i_1\cdots i_p;}{}^{j_1\cdots j_q}$, then
\begin{align}
\left((\nabla^{\mathrm H}T)_{i i_1\cdots i_p;}\right)^{
j_1\cdots j_q}
={}&H_iT_{i_1\cdots i_p;}{}^{j_1\cdots j_q}
+\omega_i^{\mathbf{F}}T_{i_1\cdots i_p;}{}^{j_1\cdots j_q}
-T_{i_1\cdots i_p;}{}^{j_1\cdots j_q}\omega_i^{\mathbf{E}}
\nonumber\\
&-\sum_{a=1}^p\Gamma^k_{i i_a}
T_{i_1\cdots k\cdots i_p;}{}^{j_1\cdots j_q}
+\sum_{b=1}^q\Gamma^{j_b}_{ik}
T_{i_1\cdots i_p;}{}^{j_1\cdots k\cdots j_q},
\label{eq:derivada-horizontal-tensor-simbolo}\\
\left((\nabla^{\mathrm V}T)_{i_1\cdots i_p;}\right)^{
j j_1\cdots j_q}
={}&D_{\xi_j}T_{i_1\cdots i_p;}{}^{j_1\cdots j_q}.
\label{eq:derivada-vertical-tensor-simbolo}
\end{align}
In particular, $\nabla^{\mathrm H}a$ takes values in $\pi^*T^*M\otimes\pi^*\operatorname{Hom}(\mathbf{E},\mathbf{F})$, whereas $\nabla^{\mathrm V}a$ takes values in $\pi^*TM\otimes\pi^*\operatorname{Hom}(\mathbf{E},\mathbf{F})$. The subscript on $\nabla_j^{\mathrm V}$ in the preceding directional formula labels the vertical direction $\mathbf{d}x^j$; it is not a covariant index of the vertical tensor. Define \[
(\nabla^{\mathrm H})^p(\nabla^{\mathrm V})^qa
\in\Gamma(\mathscr T^{p,q})
\] by applying first the $q$ vertical derivatives and then the $p$ horizontal derivatives, always using the preceding tensorial recurrence.

To record the order of differentiation without assuming commutativity, let $\mathfrak W_p:=\{1,\ldots,n\}^p$, with $\mathfrak W_0:=\{\varnothing\}$. For words $I=(i_1,\ldots,i_p)\in\mathfrak W_p$ and $J=(j_1,\ldots,j_q)\in\mathfrak W_q$, set
\begin{equation}
\nabla^{\mathrm H}_I
:=\nabla^{\mathrm H}_{i_1}\cdots\nabla^{\mathrm H}_{i_p},
\qquad
\nabla^{\mathrm V}_J
:=\nabla^{\mathrm V}_{j_1}\cdots\nabla^{\mathrm V}_{j_q}.
\label{eq:palabras-derivadas-geometricas}
\end{equation}
These expressions are deliberately directional compositions, not components of a tensorial derivative. The word preserves the order; in particular, $I$ is not replaced by a multi-index, nor are letters commuted. Only to count the vertical derivatives do we write $\beta(J):=e_{j_1}+\cdots+e_{j_q}\in\mathbb N_0^n$.

\begin{proposition}[Horizontal--vertical commutators]
\label{prop:conmutadores-horizontal-vertical-simbolos}
In a coordinate chart and the preceding local frames, the directional compositions satisfy
\begin{align}
[\nabla_j^{\mathrm V},\nabla_\ell^{\mathrm V}]a&=0,
\label{eq:conmutador-VV-simbolo}\\
[\nabla_i^{\mathrm H},\nabla_j^{\mathrm V}]a
&=-\Gamma^j_{ik}\nabla_k^{\mathrm V}a,
\label{eq:conmutador-HV-simbolo}\\
[\nabla_i^{\mathrm H},\nabla_j^{\mathrm H}]a
&=\mathbf{R}^{\mathbf{F}}(\boldsymbol{\partial}_i,\boldsymbol{\partial}_j)a
-a\mathbf{R}^{\mathbf{E}}(\boldsymbol{\partial}_i,\boldsymbol{\partial}_j)
+R_{ij\ell}{}^k\xi_k\nabla_\ell^{\mathrm V}a.
\label{eq:conmutador-HH-simbolo}
\end{align}
In particular, the order of a word can be changed only at the cost of the curvature and vertical-order terms displayed here.
\end{proposition}

\begin{proof}
The first identity follows from commutation of the derivatives $D_{\xi_j}$. Since $H_i=D_{x^i}+\Gamma^k_{i\ell}\xi_kD_{\xi_\ell}$,
\[
[H_i,D_{\xi_j}]=-\Gamma^j_{i\ell}D_{\xi_\ell},
\qquad
[H_i,H_j]=R_{ij\ell}{}^k\xi_kD_{\xi_\ell}.
\]
Adding the connection matrices gives curvature $\mathbf{R}^{\mathbf{F}}(\boldsymbol{\partial}_i,\boldsymbol{\partial}_j)a-a\mathbf{R}^{\mathbf{E}}(\boldsymbol{\partial}_i,\boldsymbol{\partial}_j)$ on $\operatorname{Hom}(\mathbf{E},\mathbf{F})$. These three identities yield \eqref{eq:conmutador-VV-simbolo}--\eqref{eq:conmutador-HH-simbolo}.
\end{proof}

The following result follows the triangular scheme of Lemma~\ref{lem:local-expression-higher-order}; here we must also record the vertical derivatives created by \eqref{eq:levantamiento-horizontal-cotangente-local}.

\begin{lemma}[Local expression for iterated geometric derivatives]
\label{lem:derivadas-geometricas-simbolos-locales}
Let $(M,\mathbf{g})$ be a smooth Riemannian manifold of dimension $n$ without boundary, not necessarily complete or compact, and let $\pi\colon T^*M\to M$. Let $\mathbf{E},\mathbf{F}\to M$ be complex vector bundles of ranks $r,s$ equipped with bundle metrics and compatible connections. The horizontal decomposition of $T(T^*M)$ is determined by the Levi--Civita connection, and the derivatives on $\pi^*\operatorname{Hom}(\mathbf{E},\mathbf{F})$ are those in \eqref{eq:definicion-derivadas-horizontal-vertical-cotangente}. Fix a chart $(U,x^1,\ldots,x^n)$ and local frames of $\mathbf{E}$ and $\mathbf{F}$, and let
\[
\mathbf{a}\in\Gamma\!\left(
\pi^*\operatorname{Hom}(\mathbf{E},\mathbf{F})|_{T^*U}
\right).
\]
Regard the coefficients as elements of $\operatorname{End}_{\mathbb C}(M_{s\times r}(\mathbb C))$. For $I\in\mathfrak W_p$ and $J\in\mathfrak W_q$, there exist coefficients $C^{I,J}_{\mu,\nu}(x,\xi)$ indexed by the explicit finite set
\begin{equation}
\mathcal A_{I,J}
:=\bigl\{(\mu,\beta(J)+\gamma)\in\mathbb N_0^n\times\mathbb N_0^n
\bigm| |\mu|+|\gamma|\leq p\bigr\},
\label{eq:indices-derivadas-geometricas-locales}
\end{equation}
such that
\begin{equation}
\nabla^{\mathrm H}_I\nabla^{\mathrm V}_J \mathbf{a}
=\sum_{(\mu,\nu)\in\mathcal A_{I,J}}
C^{I,J}_{\mu,\nu}(x,\xi)
\bigl(D_x^\mu D_\xi^\nu \mathbf{a}(x,\xi)\bigr).
\label{eq:formula-derivadas-geometricas-locales}
\end{equation}
If $\nu=\beta(J)+\gamma$, the coefficient $C^{I,J}_{\mu,\nu}$ is polynomial in $\xi$ of degree at most $|\gamma|$. For $p\geq1$, its coefficients are finite sums of products and compositions formed from $D_x^\delta\Gamma^k_{ij}$, $D_x^\delta\omega_i^{\mathbf{E}}$, and $D_x^\delta\omega_i^{\mathbf{F}}$, where $|\delta|\leq p-1$; for $p=0$ none of these data occur.

More precisely, let $L_B(T):=BT$, $R_D(T):=TD$, and \[
\mathcal H_i C
:=D_{x^i}C+\sum_{j,k=1}^n
\Gamma^k_{ij}\xi_kD_{\xi_j}C
\] for a coefficient taking values in $\operatorname{End}_{\mathbb C}(M_{s\times r}(\mathbb C))$. With the convention that a coefficient whose index does not belong to $\mathcal A_{I,J}$ is zero, the coefficients are determined by
\begin{equation}
C^{\varnothing,J}_{0,\beta(J)}
=\operatorname{Id}_{M_{s\times r}(\mathbb C)}
\label{eq:recurrencia-base-derivadas-geometricas}
\end{equation}
and, if $I^+=(i,i_1,\ldots,i_p)$,
\begin{equation}
C^{I^+,J}_{\mu,\nu}
=\mathcal H_iC^{I,J}_{\mu,\nu}
+L_{\omega_i^{\mathbf{F}}}\circ C^{I,J}_{\mu,\nu}
-R_{\omega_i^{\mathbf{E}}}\circ C^{I,J}_{\mu,\nu}+\mathbf{1}_{\{\mu_i\geq1\}}C^{I,J}_{\mu-e_i,\nu}
+\sum_{j,k=1}^n
\mathbf{1}_{\{\nu_j\geq1\}}
\Gamma^k_{ij}\xi_k C^{I,J}_{\mu,\nu-e_j}.
\label{eq:recurrencia-derivadas-geometricas}
\end{equation}
\end{lemma}

\begin{proof}
If $|I|+|J|=1$, the two possibilities are, by \eqref{eq:derivadas-geometricas-cotangente-locales},
\[
\nabla^{\mathrm V}_j \mathbf{a}=D_{\xi_j}\mathbf{a}
\]
and
\[
\nabla^{\mathrm H}_i \mathbf{a}
=D_{x^i}\mathbf{a}+\sum_{j,k=1}^{n}\Gamma^k_{ij}\xi_kD_{\xi_j}\mathbf{a}
+L_{\omega_i^{\mathbf{F}}}(\mathbf{a})-R_{\omega_i^{\mathbf{E}}}(\mathbf{a}).
\]
In the first case the only coefficient is the identity with $(\mu,\nu)=(0,e_j)$. In the second, the coefficients are the identity for $(e_i,0)$, the endomorphism $\displaystyle (\displaystyle\sum_{k=1}^{n}\Gamma^k_{ij}\xi_k)
\operatorname{Id}_{M_{s\times r}(\mathbb C)}$ for $(0,e_j)$, and $L_{\omega_i^{\mathbf{F}}}-R_{\omega_i^{\mathbf{E}}}$ for $(0,0)$. Their respective degrees are $0$, $1$, and $0$, and no derivatives of the connection data appear.

Begin with the vertical induction. For $q=1$, we have already obtained $\nabla^{\mathrm V}_j \mathbf{a}=D_{\xi_j}\mathbf{a}$. If $\nabla^{\mathrm V}_J \mathbf{a}=D_\xi^{\beta(J)}\mathbf{a}$ and $J^+=(j,j_1,\ldots,j_q)$, the first identity in \eqref{eq:derivadas-geometricas-cotangente-locales} gives
\[
\nabla^{\mathrm V}_{J^+}\mathbf{a}
=D_{\xi_j}\bigl(D_\xi^{\beta(J)}\mathbf{a}\bigr)
=D_\xi^{\beta(J)+e_j}\mathbf{a}.
\]
No additional coefficient appears because the linear coordinates on each fiber have trivial vertical connection. This proves \eqref{eq:recurrencia-base-derivadas-geometricas} for every word $J$. If a vertical letter is added to $J$ while keeping a horizontal word $I$, the convention \eqref{eq:palabras-derivadas-geometricas} places that differentiation before the horizontal ones: replace the initial datum $C^{\varnothing,J}_{0,\beta(J)}$ by $C^{\varnothing,J^+}_{0,\beta(J)+e_j}$, then apply, in the same order, the horizontal recurrences corresponding to the letters of $I$. This gives the vertical step for any $I$, with differentiation $\nabla^{\mathrm V}_j$ applied before the horizontal derivatives.

Now fix $J$ and assume \eqref{eq:formula-derivadas-geometricas-locales} for a word $I$ of length $p$. For brevity, let $A_{\mu,\nu}:=D_x^\mu D_\xi^\nu \mathbf{a}$. Applying an additional horizontal derivative to each summand and using the Leibniz rule gives the terms
\begin{align*}
\nabla_i^{\mathrm H}
\bigl(C^{I,J}_{\mu,\nu}A_{\mu,\nu}\bigr)
={}&(D_{x^i}C^{I,J}_{\mu,\nu})A_{\mu,\nu}
+\sum_{j,k=1}^{n}\Gamma^k_{ij}\xi_k
(D_{\xi_j}C^{I,J}_{\mu,\nu})A_{\mu,\nu}\\
&+C^{I,J}_{\mu,\nu}D_x^{\mu+e_i}D_\xi^\nu \mathbf{a}
+\sum_{j,k=1}^{n}\Gamma^k_{ij}\xi_k
C^{I,J}_{\mu,\nu}D_x^\mu D_\xi^{\nu+e_j}\mathbf{a}\\
&+L_{\omega_i^{\mathbf{F}}}\bigl(C^{I,J}_{\mu,\nu}A_{\mu,\nu}\bigr)
-R_{\omega_i^{\mathbf{E}}}\bigl(C^{I,J}_{\mu,\nu}A_{\mu,\nu}\bigr).
\end{align*}
The first two expressions differentiate the coefficient, the third differentiates $\mathbf{a}$ spatially, the fourth comes from the vertical part of the horizontal lift, and the last two are left and right connection multiplication, respectively. Regrouping by the derivative $D_x^\mu D_\xi^\nu \mathbf{a}$ gives exactly \eqref{eq:recurrencia-derivadas-geometricas}.

Let us verify the four inductive properties. The first three terms of the recurrence preserve $(\mu,\nu)$; the fourth replaces $\mu$ by $\mu+e_i$; the last replaces $\nu=\beta(J)+\gamma$ by $\nu+e_j$. In all three cases, $|\mu|+|\gamma|\leq p+1$. Differentiation in $x$ preserves the degree in $\xi$; the operator $\xi_kD_{\xi_j}$ does not increase it either; connection multiplication preserves it; and the term creating $D_{\xi_j}$ multiplies by a function linear in $\xi$, so the new degree is at most $|\gamma|+1$. Differentiating an earlier coefficient in $x$ may introduce derivatives of order $p$ of $\Gamma$, $\omega^{\mathbf{E}}$, or $\omega^{\mathbf{F}}$; the other terms introduce only undifferentiated data. This is precisely the bound $(p+1)-1$. Finally, \eqref{eq:indices-derivadas-geometricas-locales} is finite, and each sum in \eqref{eq:recurrencia-derivadas-geometricas} contains finitely many terms. This proves the step $p\longmapsto p+1$ and hence the lemma. The inductions $q\longmapsto q+1$ for the vertical datum and $p\longmapsto p+1$ for the horizontal letters cover, in particular, each total level $p+q=N+1$ from earlier levels.
\end{proof}

The preceding relation is triangular and can be inverted without losing control of the indices.

\begin{lemma}[Inverse local triangular formula]
\label{lem:derivadas-parciales-en-derivadas-geometricas}
Let $(M,\mathbf{g})$ be a smooth Riemannian manifold of dimension $n$ without boundary, not necessarily complete or compact. Let $\mathbf{E},\mathbf{F}\to M$ be complex vector bundles of finite rank with bundle metrics and compatible connections; use the Levi--Civita connection for the horizontal lift. Fix a chart $U$, frames of $\mathbf{E}$ and $\mathbf{F}$ over $U$, and $\mathbf{a}\in\Gamma(\pi^*\operatorname{Hom}(\mathbf{E},\mathbf{F})|_{T^*U})$. For each $\alpha,\beta\in\mathbb N_0^n$, there exist coefficients $\widetilde C^{\alpha,\beta}_{I,J}(x,\xi)$ such that
\begin{equation}
D_x^\alpha D_\xi^\beta \mathbf{a}
=\sum_{(I,J)\in\widetilde{\mathcal A}_{\alpha,\beta}}
\widetilde C^{\alpha,\beta}_{I,J}(x,\xi)
\bigl(\nabla_I^{\mathrm H}\nabla_J^{\mathrm V}\mathbf{a}\bigr),
\label{eq:formula-inversa-derivadas-geometricas}
\end{equation}
where the admissible finite set is
\begin{equation}
\widetilde{\mathcal A}_{\alpha,\beta}
:=\bigcup_{h=0}^{|\alpha|}
\bigcup_{q=|\beta|}^{|\beta|+|\alpha|-h}
\mathfrak W_h\times\mathfrak W_q.
\label{eq:indices-formula-inversa-derivadas-geometricas}
\end{equation}
For $I\in\mathfrak W_h$ and $J\in\mathfrak W_q$, the coefficient can be chosen polynomial in $\xi$ of degree at most $q-|\beta|$ and formed from derivatives in $x$ of $\Gamma$, $\omega^{\mathbf{E}}$, and $\omega^{\mathbf{F}}$ of order at most $|\alpha|-1$ when $|\alpha|\geq1$.
\end{lemma}

\begin{proof}
Order partial derivatives by their \emph{horizontal complexity} $|\alpha|$; the number of vertical derivatives does not enter the first component of this ordering. If $|\alpha|=0$, choose a word $J_\beta$ containing each index $j$ exactly $\beta_j$ times. The vertical induction in the preceding proof gives $D_\xi^\beta \mathbf{a}=\nabla^{\mathrm V}_{J_\beta}\mathbf{a}$, which is \eqref{eq:formula-inversa-derivadas-geometricas}.

Suppose the inverse formula has been constructed for every derivative of horizontal complexity less than $p$ and every vertical multi-index. Let $|\alpha|=p$ and choose a word $I_\alpha$ containing $i$ exactly $\alpha_i$ times. In \eqref{eq:formula-derivadas-geometricas-locales}, applied to $(I_\alpha,J_\beta)$, there is a unique term with $|\mu|=p$: to obtain it, each of the $p$ horizontal derivatives must choose $D_{x^i}$ acting on $\mathbf{a}$. Its coefficient is the identity, its index is $(\alpha,\beta)$, and all other terms have $|\mu|\leq p-1$. Thus
\begin{equation}
D_x^\alpha D_\xi^\beta \mathbf{a}
=\nabla^{\mathrm H}_{I_\alpha}
\nabla^{\mathrm V}_{J_\beta}\mathbf{a}
-\sum_{(\mu,\nu)\in
\mathcal A_{I_\alpha,J_\beta}\setminus\{(\alpha,\beta)\}}
C^{I_\alpha,J_\beta}_{\mu,\nu}
D_x^\mu D_\xi^\nu \mathbf{a}.
\label{eq:paso-triangular-inverso}
\end{equation}
In each summand substitute the induction hypothesis corresponding to $(\mu,\nu)$. After summing coefficients with the same words, this defines the exact recurrence
\begin{equation}
\widetilde C^{\alpha,\beta}_{I,J}
=\mathbf{1}_{\{(I,J)=(I_\alpha,J_\beta)\}}
\operatorname{Id}_{M_{s\times r}(\mathbb C)}
-\sum_{(\mu,\nu)\in\mathcal A_{I_\alpha,J_\beta}\setminus\{(\alpha,\beta)\}}
C^{I_\alpha,J_\beta}_{\mu,\nu}
\circ\widetilde C^{\mu,\nu}_{I,J},
\label{eq:recurrencia-formula-inversa-derivadas-geometricas}
\end{equation}
where the sum is restricted to $\mathcal A_{I_\alpha,J_\beta}$ and an undefined coefficient is interpreted as zero.

It remains to check the bounds. Write $\nu=\beta+\gamma$; the preceding lemma gives $|\mu|+|\gamma|\leq p$. If a word produced by substituting $D_x^\mu D_\xi^\nu \mathbf{a}$ has lengths $h$ and $q$, the induction hypothesis gives $h\leq|\mu|$ and
\[
|\nu|\leq q\leq|\nu|+|\mu|-h
\leq|\beta|+p-h.
\]
This is exactly \eqref{eq:indices-formula-inversa-derivadas-geometricas}. The degree of the product of coefficients is at most $|\gamma|+(q-|\nu|)=q-|\beta|$. Its connection jets have order at most $p-1$, since the first factor has this bound and the inductive coefficients use at most $|\mu|-1$. The sums are finite. This proves the step from horizontal complexity $p-1$ to $p$.
\end{proof}

\begin{lemma}[Tensorial derivatives and ordered words]
\label{lem:simbolos-derivadas-tensoriales-palabras}
In a chart and local frames, each component of $(\nabla^{\mathrm H})^p(\nabla^{\mathrm V})^qa$ is the directional word with the same indices plus a finite sum of words with fewer than $p$ horizontal letters. The coefficients are universal polynomials in the Christoffel symbols and their derivatives of order at most $p-1$. Conversely, each word of lengths $(p,q)$ is the corresponding component of the iterated tensor plus a sum of the same type. On each compact set, the two families of components are therefore equivalent without loss of order.
\end{lemma}

\begin{proof}
For $(p,q)=(1,0)$ or $(0,1)$ there is no correction. The first nontrivial identities are
\[
\left((\nabla^{\mathrm H})^2a\right)_{ij}
=\nabla_i^{\mathrm H}(\nabla_j^{\mathrm H}a)
-\Gamma^k_{ij}\nabla_k^{\mathrm H}a
\]
and
\[
\left(\nabla^{\mathrm H}\nabla^{\mathrm V}a\right)_{i;}{}^j
=\nabla_i^{\mathrm H}(\nabla_j^{\mathrm V}a)
+\Gamma^j_{ik}\nabla_k^{\mathrm V}a.
\]
Formulas \eqref{eq:derivada-horizontal-tensor-simbolo}--%
\eqref{eq:derivada-vertical-tensor-simbolo} show inductively that adding an index leaves as leading term the word obtained by adding the same letter, and that all corrections have smaller horizontal length. The triangular matrix has identity diagonal; solving inductively gives the converse assertion. Uniform bounds on compact sets follow from bounds on the Christoffel jets.
\end{proof}

\begin{definition}[Full geometric symbol]
\label{def:simbolo-geometrico-completo}
\index{geometric symbol@geometric symbol}
Let $(M,\mathbf{g})$ be a smooth Riemannian manifold without boundary, not necessarily complete or compact, let $\pi\colon T^*M\to M$, and let $\mathbf{E},\mathbf{F}\to M$ be complex vector bundles of respective finite ranks $r_{\mathbf{E}},r_{\mathbf{F}}$, with bundle metrics $\mathbf{h}_{\mathbf{E}},\mathbf{h}_{\mathbf{F}}$ and compatible connections. On $T^*M$ use the Levi--Civita horizontal lift, and on $\pi^*\operatorname{Hom}(\mathbf{E},\mathbf{F})$ the pullback connection induced by the connections on $\mathbf{E}$ and $\mathbf{F}$. A \textbf{geometric symbol of order $m$} is a section \[
\mathbf{a}\in \Gamma\bigl(\pi^*\operatorname{Hom}(\mathbf{E},\mathbf{F})\bigr)
\] such that, for every compact set $K\Subset M$ and all $p,q\in\mathbb N_0$, there exists $C_{K,p,q}>0$ with
\begin{equation}
\bigl\|(\nabla^{\mathrm H})^p(\nabla^{\mathrm V})^q
\mathbf{a}(x,\xi)\bigr\|_{\mathbf{g},\mathbf{h}_{\mathbf{E}},\mathbf{h}_{\mathbf{F}}}
\leq
C_{K,p,q}\langle\xi\rangle_{\mathbf{g}}^{m-q},
\qquad x\in K,\quad \xi\in T_x^*M.
\label{eq:estimaciones-simbolo-geometrico}
\end{equation}
The class is denoted by $S^m_{\mathbf{g},\nabla}(T^*M;\operatorname{Hom}(\mathbf{E},\mathbf{F}))$.

In \eqref{eq:estimaciones-simbolo-geometrico}, the notation $(\nabla^{\mathrm H})^p(\nabla^{\mathrm V})^q \mathbf{a}$ denotes the tensor in $\mathscr T^{p,q}$ defined by the preceding intrinsic recurrence, and the norm is induced by $\mathbf{g}$ and the bundle metrics. The words \eqref{eq:palabras-derivadas-geometricas} are retained only as auxiliary compositions for the local inductions.

The symbol is \textbf{classical} if there exist sections $\mathbf{a}_{m-j}$ on $T^*M\setminus0_M$, homogeneous of degree $m-j$, such that, for a radial cutoff $\vartheta$ vanishing near $0_M$ and equal to one at high frequencies,
\[
\mathbf{a}-\sum_{j=0}^{N-1}\vartheta \mathbf{a}_{m-j}
\in S^{m-N}_{\mathbf{g},\nabla}(T^*M;\operatorname{Hom}(\mathbf{E},\mathbf{F}))
\qquad(N\geq1).
\]
We then write $\mathbf{a}\in S^m_{\mathrm{cl},\mathbf{g},\nabla}
(T^*M;\operatorname{Hom}(\mathbf{E},\mathbf{F}))$.
\end{definition}

The comparison with local classes is quantitative. Fix a matrix norm $\|\cdot\|_0$ in the chosen frames, and denote by $\|\cdot\|_{\operatorname{op},0}$ the induced norm on the corresponding endomorphism spaces. For a tensor $T\in\mathscr T^{p,q}$, let $\|T\|_{0,\mathrm{comp}}$ denote the Euclidean product norm of the matrix norms of its components in the fixed chart and frames; for $(p,q)=(0,0)$ it agrees with $\|\cdot\|_0$. For $K\Subset U$, $N\in\mathbb N_0$, and $m\in\mathbb R$, define
\begin{align}
p^{\mathrm{loc}}_{K,N,m}(a)
&:=\max_{|\alpha|+|\beta|\leq N}
\sup_{x\in K,\,\xi\in\mathbb R^n}
\langle\xi\rangle^{-m+|\beta|}
\|D_x^\alpha D_\xi^\beta a(x,\xi)\|_0,
\label{eq:seminorma-local-simbolo-geometrico}\\
p^{\mathrm{geom}}_{K,N,m}(a)
&:=\max_{p+q\leq N}
\sup_{x\in K,\,\xi\in T_x^*M}
\langle\xi\rangle_{\mathbf{g}}^{-m+q}
\|(\nabla^{\mathrm H})^p(\nabla^{\mathrm V})^qa(x,\xi)\|_{\mathbf{g},\mathbf{h}_{\mathbf{E}},\mathbf{h}_{\mathbf{F}}},
\label{eq:seminorma-geometrica-simbolo}\\
p^{\mathrm{word}}_{K,N,m}(a)
&:=\max_{\substack{p+q\leq N\\
I\in\mathfrak W_p,\ J\in\mathfrak W_q}}
\sup_{x\in K,\,\xi\in T_x^*M}
\langle\xi\rangle_{\mathbf{g}}^{-m+q}
\|\nabla_I^{\mathrm H}\nabla_J^{\mathrm V}a(x,\xi)\|_{\mathbf{g},\mathbf{h}_{\mathbf{E}},\mathbf{h}_{\mathbf{F}}}.
\label{eq:seminorma-palabras-simbolo}
\end{align}
The tensor norms are induced by $\mathbf{g}$ and the bundle metrics.

\begin{proposition}[Quantitative equivalence of symbol seminorms]
\label{prop:equivalencia-seminormas-simbolos-geometricos-locales}
Let $(M,\mathbf{g})$ be a smooth Riemannian manifold without boundary, not necessarily complete or compact, and let $\mathbf{E},\mathbf{F}\to M$ be complex vector bundles of finite rank with bundle metrics $\mathbf{h}_{\mathbf{E}},\mathbf{h}_{\mathbf{F}}$ and compatible connections. Fix a chart $(U,x)$, local frames of $\mathbf{E}$ and $\mathbf{F}$, and a symbol $\mathbf{a}\in\Gamma(\pi^*\operatorname{Hom}(\mathbf{E},\mathbf{F})|_{T^*U})$. For each $K\Subset U$, $N\in\mathbb N_0$, and $m\in\mathbb R$, there exists $C_{K,N,m}\geq1$ such that
\begin{equation}
\begin{aligned}
C_{K,N,m}^{-1}p^{\mathrm{loc}}_{K,N,m}(\mathbf{a})
&\leq p^{\mathrm{word}}_{K,N,m}(\mathbf{a})
\leq C_{K,N,m}p^{\mathrm{loc}}_{K,N,m}(\mathbf{a}),\\
C_{K,N,m}^{-1}p^{\mathrm{word}}_{K,N,m}(\mathbf{a})
&\leq p^{\mathrm{geom}}_{K,N,m}(\mathbf{a})
\leq C_{K,N,m}p^{\mathrm{word}}_{K,N,m}(\mathbf{a}).
\end{aligned}
\label{eq:equivalencia-cuantitativa-seminormas-simbolos}
\end{equation}
In particular, the same order $N$ suffices in both inequalities: no derivatives are lost. The constants depend on $K,N,m$, the chart and frames, comparison of the Euclidean and Riemannian norms, the bundle metrics, and bounds on $K$ for derivatives of $\Gamma$, $\omega^{\mathbf{E}}$, and $\omega^{\mathbf{F}}$ of order at most $N-1$; when $N=0$ they do not depend on these jets.
\end{proposition}

\begin{proof}
Compactness of $K$ provides constants $c_K,C_K,c'_K,C'_K>0$ with
\[
c_K\langle\xi\rangle
\leq\langle\xi\rangle_{\mathbf{g}}
\leq C_K\langle\xi\rangle,
\qquad
c'_K\|T\|_{0,\mathrm{comp}}
\leq\|T\|_{\mathbf{g},\mathbf{h}_{\mathbf{E}},\mathbf{h}_{\mathbf{F}}}\leq
C'_K\|T\|_{0,\mathrm{comp}}.
\]
The polynomial and jet bounds in Lemma~\ref{lem:derivadas-geometricas-simbolos-locales} imply that, for $p+q\leq N$, there exists $B_{K,N}>0$ such that
\[
\|C^{I,J}_{\mu,\beta(J)+\gamma}(x,\xi)\|_{\operatorname{op},0}
\leq B_{K,N}\langle\xi\rangle^{|\gamma|},
\qquad x\in K,\quad \xi\in\mathbb R^n.
\]
Since $|\mu|+|\beta(J)+\gamma|\leq p+q\leq N$, each partial derivative of \eqref{eq:formula-derivadas-geometricas-locales} is included in $p^{\mathrm{loc}}_{K,N,m}$. Moreover,
\[
\langle\xi\rangle^{|\gamma|}
\|D_x^\mu D_\xi^{\beta(J)+\gamma}\mathbf{a}\|_0
\leq p^{\mathrm{loc}}_{K,N,m}(\mathbf{a})
\langle\xi\rangle^{m-q}.
\]
Summing over the finitely many indices in \eqref{eq:indices-derivadas-geometricas-locales} and using both norm comparisons proves the bound of $p^{\mathrm{word}}$ by $p^{\mathrm{loc}}$ in \eqref{eq:equivalencia-cuantitativa-seminormas-simbolos}.

For the reverse direction, Lemma~\ref{lem:derivadas-parciales-en-derivadas-geometricas} gives a constant $\widetilde B_{K,N}>0$ such that, if $(I,J)\in\mathfrak W_h\times\mathfrak W_q$ contributes to $D_x^\alpha D_\xi^\beta \mathbf{a}$ with $|\alpha|+|\beta|\leq N$, then
\[
\|\widetilde C^{\alpha,\beta}_{I,J}(x,\xi)\|_{\operatorname{op},0}
\leq\widetilde B_{K,N}\langle\xi\rangle^{q-|\beta|},
\qquad
h+q\leq|\alpha|+|\beta|\leq N.
\]
By definition of the word seminorm,
\[
\langle\xi\rangle^{q-|\beta|}
\|\nabla_I^{\mathrm H}\nabla_J^{\mathrm V}\mathbf{a}\|_{\mathbf{g},\mathbf{h}_{\mathbf{E}},\mathbf{h}_{\mathbf{F}}}
\leq C_{K,N,m}'p^{\mathrm{word}}_{K,N,m}(\mathbf{a})
\langle\xi\rangle^{m-|\beta|}.
\]
The finite sum in \eqref{eq:formula-inversa-derivadas-geometricas} bounds $p^{\mathrm{loc}}$ by $p^{\mathrm{word}}$. Recurrences \eqref{eq:recurrencia-derivadas-geometricas} and \eqref{eq:recurrencia-formula-inversa-derivadas-geometricas} show exactly the asserted dependence of the constants. Finally, Lemma~\ref{lem:simbolos-derivadas-tensoriales-palabras}, applied for $p+q\leq N$, compares $p^{\mathrm{geom}}$ with $p^{\mathrm{word}}$ at the same orders. Increasing the constant once gives all four stated inequalities.
\end{proof}

Proposition~\ref{prop:equivalencia-seminormas-simbolos-geometricos-locales} shows that, in every relatively compact chart $(U,\kappa)$ and local bundle frames, $S^m_{\mathbf{g},\nabla}(T^*M;\operatorname{Hom}(\mathbf{E},\mathbf{F}))$ coincides with $S^m_{1,0}(\kappa(U)\times\mathbb R^n;
\operatorname{Hom}(\mathbb C^{r_{\mathbf{E}}},\mathbb C^{r_{\mathbf{F}}}))$, and that both families of seminorms generate the same Fréchet topology. It also shows that replacing the metric, connections, frames, or charts by other smooth data does not change the local class: apply \eqref{eq:equivalencia-cuantitativa-seminormas-simbolos} to each choice and compare through the local family. On a noncompact base, all constants depend on $K$, so this equivalence is uniform only on compact sets.

\begin{definition}[Riemannian quantization]
\label{def:cuantizacion-riemanniana-pseudo}
\index{quantization@quantization!Riemannian pseudodifferential}
Let $(M,\mathbf{g})$ be a smooth Riemannian manifold of dimension $n$ without boundary, not necessarily complete or compact. Let $\mathbf{E},\mathbf{F}\to M$ be complex vector bundles of finite rank with bundle metrics and compatible connections; let $\chi\in C^\infty(M\times M)$ be the properly supported cutoff fixed at the beginning of this section, and use parallel transport in $\mathbf{E}$ in the input variable. For $a\in S^m_{\mathbf{g},\nabla}(T^*M;\operatorname{Hom}(\mathbf{E},\mathbf{F}))$ and $\mathbf{u}\in\Gamma_c(M,\mathbf{E})$, define
\[
\begin{split}
\operatorname{Op}_{\mathbf{g},\nabla,\chi}(a)\mathbf{u}(x)
:={}&
(2\pi)^{-n}\operatorname{Os}\!\int_{T_x^*M}\int_{T_xM}
e^{-i\langle\xi,v\rangle}
\chi(x,\exp_xv)\,a(x,\xi)\\
&\hspace{35mm}\cdot
\tau^{\mathbf{E}}_{x\leftarrow\exp_xv}\mathbf{u}(\exp_xv)\,dv_x\,d\xi_x.
\end{split}
\]
\end{definition}

The formula can be written using Riemannian measure. If \[
\exp_x^*(d\lambda_{\mathbf{g}})(v)=j_x(v)\,dv_x,
\], the integrand over $M$ contains the factor $j_x(\ell(x,y))^{-1}d\lambda_{\mathbf{g}}(y)$. Thus Definition~\ref{def:cuantizacion-riemanniana-pseudo} is compatible with the global convention $d\lambda_{\mathbf{g}}$: the factor $j_x(\ell(x,y))^{-1}$ is the Jacobian relating the two measures.

This construction is the pseudodifferential quantization associated with the linearization $\ell$, not the geometric quantization of Hamiltonian systems. It follows the intrinsic idea of Bokobza--Haggiag \cite[Chapters~III--VI]{BokobzaHaggiag1969}: a linearization, a cutoff, and a transport turn a cotangent symbol into a kernel near the diagonal.

\begin{theorem}[Equivalence of local and Riemannian quantization]
\label{teo:equiv-local-geometrica-pseudo}
Let $(M,\mathbf{g})$ be a smooth Riemannian manifold without boundary, not necessarily complete or compact. Let $\mathbf{E},\mathbf{F}\to M$ be complex vector bundles of finite rank with bundle metrics and compatible connections, and let $\chi$ be a properly supported cutoff as in Definition~\ref{def:cuantizacion-riemanniana-pseudo}. The following two assertions hold.
\begin{enumerate}[label=(\roman*)]
\item If $a\in S^m_{\mathrm{cl},\mathbf{g},\nabla}
(T^*M;\operatorname{Hom}(\mathbf{E},\mathbf{F}))$, then \[
\operatorname{Op}_{\mathbf{g},\nabla,\chi}(a)
\in\Psi^m_{\mathrm{cl}}(M;\mathbf{E},\mathbf{F})
\], and its principal symbol is the homogeneous component $a_m$.
\item If $A\in\Psi^m_{\mathrm{cl}}(M;\mathbf{E},\mathbf{F})$, there exists a classical geometric symbol $a$ such that
\[
A-\operatorname{Op}_{\mathbf{g},\nabla,\chi}(a)
\in\Psi^{-\infty}(M;\mathbf{E},\mathbf{F}).
\]
If the smooth kernel of this remainder is allowed to be added to the quantization, reconstruction is exact.
\end{enumerate}
\end{theorem}

\begin{proof}
Throughout this proof, $M,\mathbf{E},\mathbf{F}$ and the geometric data remain fixed; we abbreviate $\Psi^q:=\Psi^q(M;\mathbf{E},\mathbf{F})$, $\Psi^{-\infty}:=\Psi^{-\infty}(M;\mathbf{E},\mathbf{F})$, and $S^q_{\mathbf{g},\nabla}:=S^q_{\mathbf{g},\nabla}(T^*M;
\operatorname{Hom}(\mathbf{E},\mathbf{F}))$. For (i), localize the operator with functions supported in a chart, choosing the chart so that the segments appearing below remain inside it. Let $\mathbf{g}_0(x):=\det(\mathbf{g}(x))$, and let $T(x,y)$ be the matrix of $\tau^{\mathbf{E}}_{x\leftarrow y}$. Since $\ell(x,x)=0$, the parameter-dependent version of Hadamard's lemma~\ref{lem:hadamard-finito-dimensional} gives, componentwise, the exact identity
\begin{equation}
\ell^a(x,y)=\sum_{b=1}^nC^a_b(x,y)(y^b-x^b),
\qquad
C^a_b(x,y):=\int_0^1
\partial_{y^b}\ell^a(x,x+t(y-x))\,dt.
\label{eq:factorizacion-logaritmo-riemanniano}
\end{equation}
Since $D_y\ell(x,x)=\operatorname{Id}_{T_xM}$, we have $C(x,x)=I$. After shrinking the localization support, $C(x,y)$ is invertible, and $C$, $C^{-1}$, and every finite family of their derivatives are bounded.

In cotangent coordinates, $d\xi_x=\mathbf{g}_0(x)^{-\displaystyle\frac{1}{2}}d\xi$, whereas $d\lambda_{\mathbf{g}}(y)=\mathbf{g}_0(y)^{\displaystyle\frac{1}{2}}dy$. The formula preceding the theorem and \eqref{eq:factorizacion-logaritmo-riemanniano} turn the localized operator into
\[
(2\pi)^{-n}
\operatorname{Os}\!\iint_{\mathbb R_y^n\times\mathbb R_\xi^n}
e^{-i(y-x)\cdot C(x,y)^T\xi}
h(x,y)a(x,\xi)T(x,y)u(y)\,dy\,d\xi,
\]
where, including the two localization functions in $h$,
\[
h(x,y)=\chi(x,y)j_x(\ell(x,y))^{-1}
\left(\frac{\mathbf{g}_0(y)}{\mathbf{g}_0(x)}\right)^{\frac{1}{2}}.
\]
Now make the exact change
\begin{equation}
\eta=C(x,y)^T\xi,
\qquad
\xi=L(x,y)\eta,
\qquad
L(x,y):=C(x,y)^{-T}.
\label{eq:cambio-covariable-cuantizacion-riemanniana}
\end{equation}
Its Jacobian is $d\xi=|\det C(x,y)|^{-1}d\eta$. Consequently, the full local amplitude is
\begin{equation}
b(x,y,\eta)
=h(x,y)|\det C(x,y)|^{-1}
a\bigl(x,L(x,y)\eta\bigr)T(x,y).
\label{eq:amplitud-local-cuantizacion-riemanniana}
\end{equation}

Let us verify its estimates. The chain rules, applying derivatives in $\eta$ first and then spatial derivatives, produce a finite sum of terms
\begin{equation}
B_{\alpha,\beta,\gamma;\mu,\nu,\rho}(x,y)
\eta^\rho
(D_x^\mu D_\xi^\nu a)\bigl(x,L(x,y)\eta\bigr),
\label{eq:cadena-amplitud-geometrica-local}
\end{equation}
when computing $D_x^\alpha D_y^\beta D_\eta^\gamma[a(x,L\eta)]$. The indices in each term satisfy
\[
|\nu|=|\gamma|+|\rho|,
\qquad
|\mu|+|\rho|\leq|\alpha|+|\beta|,
\]
and the coefficients $B_{\alpha,\beta,\gamma;\mu,\nu,\rho}$ are products of derivatives of $L$. These assertions follow inductively: a spatial derivative hitting $a$ produces $D_xa$; one hitting $L\eta$ simultaneously produces a linear factor in $\eta$ and an additional derivative in $\xi$; and one hitting a coefficient merely differentiates that coefficient. This enumerates the three Leibniz possibilities and preserves both displayed relations.

On the spatial compact set under consideration, there exist $c,C>0$ such that $c\langle\eta\rangle\leq\langle L(x,y)\eta\rangle\leq
C\langle\eta\rangle$. Proposition \ref{prop:equivalencia-seminormas-simbolos-geometricos-locales} and \eqref{eq:cadena-amplitud-geometrica-local} give, for all $\alpha,\beta,\gamma$, a constant $C_{\alpha,\beta,\gamma}>0$ with
\[
\bigl\|D_x^\alpha D_y^\beta D_\eta^\gamma
a(x,L\eta)\bigr\|
\leq C_{\alpha,\beta,\gamma}
\langle\eta\rangle^{m-|\gamma|}.
\]
Indeed, the factor $\eta^\rho$ in each summand is exactly offset by the $|\rho|$ additional derivatives in $\xi$. All derivatives of $h$, $|\det C|^{-1}$, and $T$ are bounded on the compact set, so
\begin{equation}
\|D_x^\alpha D_y^\beta D_\eta^\gamma b(x,y,\eta)\|
\leq C'_{\alpha,\beta,\gamma}
\langle\eta\rangle^{m-|\gamma|}.
\label{eq:estimacion-amplitud-cuantizacion-riemanniana}
\end{equation}
This proves, derivative by derivative, that $b$ is an amplitude of order $m$. If $a$ is classical, substituting $L\eta$ into each homogeneous component preserves its degree, and the same estimate applied to each remainder proves that $b$ is classical. Theorem~\ref{teo:amplitud-a-simbolo} reduces \eqref{eq:amplitud-local-cuantizacion-riemanniana} to a local classical symbol of order $m$.

On the diagonal, $C(x,x)=L(x,x)=I$, $T(x,x)=I$, $j_x(0)=1$, and $\mathbf{g}_0(x)^{\displaystyle\frac{1}{2}}\mathbf{g}_0(x)^{-\displaystyle\frac{1}{2}}=1$ hold, and the localization cutoffs can be taken equal to one. Thus $b_m(x,x,\eta)=a_m(x,\eta)$. The zero-index term in amplitude-to-symbol reduction preserves this component; every remaining term contains at least one derivative in $\eta$ and has order at most $m-1$. This proves (i), including the statement on the principal symbol.

For (ii), we construct all terms, not just the first. Set $R_0:=A$ and $r_m:=\boldsymbol{\sigma}_m^\Psi(R_0)$. A radial cutoff, vanishing near the zero section and equal to one for large covectors, extends $r_m$ to a classical geometric symbol $a^{(0)}$ of order $m$. By (i), $\boldsymbol{\sigma}_m^\Psi(\operatorname{Op}_{\mathbf{g},\nabla,\chi}(a^{(0)}))=r_m$; hence
\[
R_1:=R_0-\operatorname{Op}_{\mathbf{g},\nabla,\chi}(a^{(0)})
\in\Psi^{m-1}_{\mathrm{cl}}.
\]

The precise induction statement is as follows: for some $N\geq1$, we have constructed $a^{(j)}\in S^{m-j}_{\mathrm{cl},\mathbf{g},\nabla}$, $0\leq j<N$, so that
\begin{equation}
R_N:=A-\operatorname{Op}_{\mathbf{g},\nabla,\chi}
\left(\sum_{j=0}^{N-1}a^{(j)}\right)
\in\Psi^{m-N}_{\mathrm{cl}}.
\label{eq:induccion-reconstruccion-simbolo-geometrico}
\end{equation}
The case $N=1$ has just been proved. Assume \eqref{eq:induccion-reconstruccion-simbolo-geometrico} and let $r_{m-N}:=\boldsymbol{\sigma}_{m-N}^\Psi(R_N)$. Extend $r_{m-N}$ by a radial cutoff to $a^{(N)}\in S^{m-N}_{\mathrm{cl},\mathbf{g},\nabla}$. Part (i) gives
\[
\boldsymbol{\sigma}_{m-N}^\Psi
\bigl(\operatorname{Op}_{\mathbf{g},\nabla,\chi}(a^{(N)})\bigr)=r_{m-N}.
\]
By exactness of the principal symbol,
\[
R_{N+1}:=R_N-\operatorname{Op}_{\mathbf{g},\nabla,\chi}(a^{(N)})
\in\Psi^{m-N-1}_{\mathrm{cl}},
\]
which is precisely \eqref{eq:induccion-reconstruccion-simbolo-geometrico} with $N+1$. This completes the induction.

Lemma~\ref{lem:suma-asintotica-simbolos}, applied in charts and combined using a locally finite partition, provides $a\in S^m_{\mathrm{cl},\mathbf{g},\nabla}$ such that, for every $N\geq1$,
\begin{equation}
a-\sum_{j=0}^{N-1}a^{(j)}
\in S^{m-N}_{\mathbf{g},\nabla}.
\label{eq:resto-suma-asintotica-simbolo-geometrico}
\end{equation}
Part (i), now at order $m-N$, turns the left-hand side of \eqref{eq:resto-suma-asintotica-simbolo-geometrico} into an operator in $\Psi^{m-N}$. Adding it to $R_N$ gives
\[
A-\operatorname{Op}_{\mathbf{g},\nabla,\chi}(a)\in\Psi^{m-N}
\qquad\text{for every }N\geq1.
\]
The intersection of these classes is $\Psi^{-\infty}$ by the characterization through smooth kernels. This proves (ii).
\end{proof}

The theorem distinguishes two levels that must not be confused: the principal component is recovered in a single step, modulo $\Psi^{m-1}(M;\mathbf{E},\mathbf{F})$; equivalence modulo $\Psi^{-\infty}(M;\mathbf{E},\mathbf{F})$ requires constructing all terms and applying asymptotic summation.

\begin{proposition}[Dependence on choices]
\label{prop:dependencia-elecciones-cuantizacion-geometrica}
Let $M$ be a finite-dimensional smooth manifold without boundary, and let $\mathbf{E},\mathbf{F}\to M$ be complex vector bundles of finite rank. Consider two geometric quantization rules obtained by choosing Riemannian metrics, bundle metrics, compatible connections, linearizations, and properly supported cutoffs as in Definition~\ref{def:cuantizacion-riemanniana-pseudo}. For a symbol of order $m$, changing the cutoff changes the quantization by an operator in $\Psi^{-\infty}(M;\mathbf{E},\mathbf{F})$. Changing the other data modifies the full symbol but preserves the class $\Psi^m(M;\mathbf{E},\mathbf{F})$, the principal symbol, and quantization modulo $\Psi^{m-1}(M;\mathbf{E},\mathbf{F})$. After successive corrections of lower-order terms, the two quantizations correspond modulo $\Psi^{-\infty}(M;\mathbf{E},\mathbf{F})$.
\end{proposition}

\begin{proof}
Throughout this proof, $M,\mathbf{E},\mathbf{F}$ remain fixed; write $\Psi^q:=\Psi^q(M;\mathbf{E},\mathbf{F})$ and $\Psi^{-\infty}:=\Psi^{-\infty}(M;\mathbf{E},\mathbf{F})$. Work in a relatively compact chart and write $v=y-x$. We use the exact factorizations in Hadamard's lemma~\ref{lem:hadamard-finito-dimensional}, with $x$ as parameter. If a smooth matrix-valued function $F(x,v)$ satisfies $F(x,0)=0$, then
\begin{equation}
F(x,v)=\sum_{j=1}^nv_jF_j(x,v),
\qquad
F_j(x,v)=\int_0^1\partial_{v_j}F(x,tv)\,dt.
\label{eq:hadamard-dependencia-cuantizacion}
\end{equation}
If, in addition, $D_vF(x,0)=0$, then
\begin{equation}
F(x,v)=\sum_{|\alpha|=2}v^\alpha F_\alpha(x,v),
\qquad
F_\alpha(x,v)=\frac{2}{\alpha!}
\int_0^1(1-t)\partial_v^\alpha F(x,tv)\,dt.
\label{eq:taylor-segundo-orden-dependencia-cuantizacion}
\end{equation}
In general, if the vertical derivatives of orders less than $L$ vanish at $v=0$, the same proof gives
\begin{equation}
F(x,v)=\sum_{|\alpha|=L}v^\alpha F_{\alpha,L}(x,v),
\qquad
F_{\alpha,L}(x,v)=\frac{L}{\alpha!}
\int_0^1(1-t)^{L-1}\partial_v^\alpha F(x,tv)\,dt.
\label{eq:taylor-orden-arbitrario-dependencia-cuantizacion}
\end{equation}
These identities also provide the required bounds. For example, for $K$ compact and $N\in\mathbb N_0$, there exists $C_{K,N,L}>0$ such that
\[
\max_{|\rho|+|\sigma|\leq N}
\sup_{x\in K,\ v\in\overline B_{\mathrm{euc}}(0,\varepsilon)}
\|D_x^\rho D_v^\sigma F_{\alpha,L}(x,v)\|
\leq C_{K,N,L}
\max_{|\rho'|+|\sigma'|\leq N+L}
\sup_{x\in K,\ v\in\overline B_{\mathrm{euc}}(0,\varepsilon)}
\|D_x^{\rho'}D_v^{\sigma'}F(x,v)\|.
\]

Each transverse monomial translates exactly into a change of symbol order. For the normalized linear phase,
\begin{equation}
v^\alpha e^{-i\eta\cdot v}
=i^{|\alpha|}D_\eta^\alpha e^{-i\eta\cdot v},
\qquad
\int_{\mathbb R_\eta^n} e^{-i\eta\cdot v}v^\alpha c(\eta)\,d\eta
=(-i)^{|\alpha|}\int_{\mathbb R_\eta^n} e^{-i\eta\cdot v}
D_\eta^\alpha c(\eta)\,d\eta.
\label{eq:factor-transversal-derivada-covariable}
\end{equation}
The second identity is understood first with a cutoff in $\eta$ and then through the oscillatory procedure already established. Thus a factor $v^\alpha$ applied to an amplitude of order $r$ produces an amplitude of order $r-|\alpha|$.

First let $\chi_0$ and $\chi_1$ be two admissible cutoffs. Their difference vanishes on a neighborhood of $v=0$, so all its vertical derivatives at $v=0$ are zero. For each $L$, formulas \eqref{eq:taylor-orden-arbitrario-dependencia-cuantizacion} and \eqref{eq:factor-transversal-derivada-covariable} express the difference kernel through an amplitude of order $m-L$. Given any prescribed order, choose $L$ larger; the difference belongs to every order and therefore has smooth kernel by the characterization of smoothing operators.

Now consider two connections on $\mathbf{E}$ and their radial transports $T_0(x,x+v)$ and $T_1(x,x+v)$. Both equal the identity at $v=0$. Define
\[
F(x,v):=T_1(x,x+v)-T_0(x,x+v).
\]
Applied to this function $F$, identity \eqref{eq:hadamard-dependencia-cuantizacion} gives
\begin{equation}
T_1(x,x+v)-T_0(x,x+v)
=\sum_{j=1}^nv_jT_j(x,v),
\qquad
T_j(x,v)=\int_0^1
\partial_{v_j}F(x,tv)\,dt.
\label{eq:diferencia-transportes-hadamard}
\end{equation}
The coefficients $T_j$ satisfy the preceding bounds for each finite family of derivatives. Integration by parts in \eqref{eq:factor-transversal-derivada-covariable} shows that the difference between the quantizations, with all other data fixed, has amplitude of order $m-1$.

The same factorization controls Jacobian factors and measures. In the coordinates of \eqref{eq:amplitud-local-cuantizacion-riemanniana}, the scalar combination corresponding to any geometric data is
\[
q(x,y)=j_x(\ell(x,y))^{-1}
\left(\frac{g_0(y)}{g_0(x)}\right)^{\frac{1}{2}}
|\det C(x,y)|^{-1},
\qquad q(x,x)=1.
\]
For a nonexponential linearization, $q$ denotes the exact factor obtained from the Jacobian of $y\longmapsto\ell(x,y)$ and the chosen densities; identity $D_y\ell(x,x)=I$ again gives $q(x,x)=1$. Thus, for two choices, $q_1-q_0$ satisfies
\eqref{eq:hadamard-dependencia-cuantizacion}
its coefficients and all their derivatives up to order $N$ are bounded by a constant $C_{K,N}$ on each compact set. Each term of the difference again has order at most $m-1$ after \eqref{eq:factor-transversal-derivada-covariable}.

It remains to compare two linearizations $\ell_0$ and $\ell_1$. Both satisfy
\[
\ell_\kappa(x,x)=0,
\qquad
D_y\ell_\kappa(x,x)=I,
\qquad \kappa\in\{0,1\}.
\]
Consequently, $G(x,v):=\ell_1(x,x+v)-\ell_0(x,x+v)$ satisfies both hypotheses of \eqref{eq:taylor-segundo-orden-dependencia-cuantizacion}, and
\begin{equation}
G(x,v)=\sum_{|\alpha|=2}v^\alpha G_\alpha(x,v),
\qquad
G_\alpha(x,v)=\frac{2}{\alpha!}\int_0^1(1-t)
\partial_v^\alpha G(x,tv)\,dt.
\label{eq:diferencia-linealizaciones-taylor}
\end{equation}
Define $\ell_t=\ell_0+tG$. Shrinking the neighborhood ensures that all $\ell_t$, $0\leq t\leq1$, are linearizations and that the matrices $C_t$ in \eqref{eq:factorizacion-logaritmo-riemanniano} are uniformly invertible. The fundamental theorem of calculus gives the exact identity
\begin{equation}
e^{-i\langle\xi,\ell_1\rangle}
-e^{-i\langle\xi,\ell_0\rangle}
=-i\int_0^1e^{-i\langle\xi,\ell_t\rangle}
\langle\xi,G\rangle\,dt
\label{eq:diferencia-fases-linealizaciones}
\end{equation}
In each integrand make the change $\eta=C_t(x,y)^T\xi$. By \eqref{eq:diferencia-linealizaciones-taylor}, the factor $\langle\xi,G\rangle$ becomes a sum of monomials $v^\alpha$, $|\alpha|=2$, multiplied by a function linear in $\eta$ and smooth coefficients uniformly bounded in $t$. Before integration by parts, the amplitude has order $m+1$; two applications of \eqref{eq:factor-transversal-derivada-covariable} bring it to order $m-1$. The Leibniz rules used in \eqref{eq:cadena-amplitud-geometrica-local} give, for all $\rho,\sigma,\gamma$, a constant $C_{K,\rho,\sigma,\gamma}>0$ independent of $t$ such that the final amplitude $d_t$ satisfies
\begin{equation}
\|D_x^\rho D_y^\sigma D_\eta^\gamma d_t(x,y,\eta)\|
\leq C_{K,\rho,\sigma,\gamma}
\langle\eta\rangle^{m-1-|\gamma|}.
\label{eq:estimacion-diferencia-linealizaciones}
\end{equation}
Integration in $t\in[0,1]$ preserves the estimate.

A change of metric splits into a change of its exponential linearization, already treated, and a change of the factors $q$; a change of connection on $\mathbf{E}$ reduces to \eqref{eq:diferencia-transportes-hadamard}. The connection on $\mathbf{F}$ enters the symbol seminorms, not the integral transporting the input section, and Proposition \ref{prop:equivalencia-seminormas-simbolos-geometricos-locales} shows that changing it preserves the symbol class. Combining \eqref{eq:diferencia-transportes-hadamard}, the factorization of $q$, and \eqref{eq:estimacion-diferencia-linealizaciones} gives, for each finite family of derivatives, the estimate for an amplitude of order $m-1$. Thus all the quantizations considered have the same principal symbol; this is a consequence of the factorizations, not a hypothesis.

Finally, let $Q_0,Q_1$ be two such quantization rules. What we have proved says that $Q_0(a)-Q_1(a)\in\Psi^{m-1}$ for every symbol of order $m$. Take the principal symbol of this difference, choose a correction of order $m-1$ for one of the two symbols, and subtract. If after $N$ corrections the difference belongs to $\Psi^{m-N}$, its principal symbol of degree $m-N$ determines the next correction, and the new difference belongs to $\Psi^{m-N-1}$. An asymptotic sum of the corrections produces two corresponding symbols whose operator difference belongs to $\Psi^{m-N}$ for every $N$, that is, to $\Psi^{-\infty}$.
\end{proof}

This separates what is canonical from what is auxiliary. The class $\Psi^m(M;\mathbf{E},\mathbf{F})$, principal symbol, ellipticity, and operator class modulo lower orders are canonical. The directional microlocal invariant will be added after it has been constructed by Fourier methods and its invariance proved. The full symbol, specific quantization, metric, connections, cutoff, and transport depend on choices. Calling the Riemannian formula \emph{intrinsic} means that all these data are geometric and global, not that the result is independent of them term by term.

On a curved manifold, replacing partial derivatives by covariant derivatives in the flat composition formula requires additional lower-order terms. These terms contain jets of the exponential map, the connection, and curvature. The particularly simple formulas in \cite[Chapter~VI]{BokobzaHaggiag1969} require a locally affine linearization; the existence of such a linearization is equivalent to a local affine structure. For the general calculus, the local formulas in Theorems~\ref{teo:composicion-simbolica-pseudo} and~\ref{teo:simbolo-adjunto-pseudo}, together with the equivalence just proved, provide expansions containing all these geometric terms.

\section{Elliptic inversion and microlocal analysis}
\label{sec:pseudo-parametricas-microlocal}

The symbolic calculus was constructed to invert elliptic operators without requiring an exact inverse. We first invert the principal component and successively correct lower-order terms. The final remainder will be smoothing; this arbitrary gain of derivatives is the source of both elliptic regularity and its microlocal version.

\begin{definition}[Parametrix]
\label{def:parametrix-pseudodiferencial}
\index{parametrix}
\glsadd{parametrix}
Let $M$ be a finite-dimensional smooth manifold without boundary, and let $\mathbf{E},\mathbf{F}\to M$ be smooth complex vector bundles of finite rank. Let $A\in\Psi^m(M;\mathbf{E},\mathbf{F})$. An operator $B\in\Psi^{-m}(M;\mathbf{F},\mathbf{E})$ is a \textbf{right parametrix} if $AB-I_{\mathbf{F}}\in\Psi^{-\infty}(M;\mathbf{F},\mathbf{F})$, a \textbf{left parametrix} if $BA-I_{\mathbf{E}}\in\Psi^{-\infty}(M;\mathbf{E},\mathbf{E})$, and a \textbf{two-sided parametrix} if both conditions hold. For the purely local calculus, replace all three classes by their variants with the subscript $\mathrm{loc}$.
\end{definition}

Ellipticity allows inversion of an operator at high frequencies but does not guarantee an algebraic or functional inverse. The distinction is essential: the first assertion produces a parametrix modulo smoothing operators; the second also requires control of the kernel and cokernel in a specific realization. We first develop the symbolic construction and then measure directional regularity using the wavefront set. All proofs use the convention $D_j=\partial_j$ and hence the factors $i^{-|\alpha|}$ in Theorem~\ref{teo:composicion-simbolica-pseudo}.

\subsection{Criteria for ellipticity}

We begin by collecting several formulations that will be used interchangeably. For a matrix $T$, let $s_{\min}(T)$ denote its smallest singular value.

\begin{proposition}[Equivalent formulations of ellipticity]
\label{prop:equivalencias-elipticidad-pseudo}
Let $U\subseteq\mathbb R^n$ be open, and let $r\in\mathbb N$ and $m\in\mathbb R$. Let
\[
a\in S^m_{\mathrm{cl}}
\bigl(U;\operatorname{End}(\mathbb C^r)\bigr),
\qquad
a\sim\sum_{j=0}^{\infty}a_{m-j}.
\]
The following conditions are equivalent, with constants that may depend on each compact set $K\Subset U$:
\begin{enumerate}[label=(\alph*)]
\item $a$ is elliptic in the sense of Definition~\ref{def:simbolo-eliptico-euclidiano};
\item there exist $c_K,R_K>0$ such that
\[
\|a(x,\xi)v\|
\geq c_K\langle\xi\rangle^m\|v\|,
\qquad
x\in K,\quad \|\xi\|\geq R_K;
\]
\item there exist $c_K,R_K>0$ such that
\[
|\det a(x,\xi)|
\geq c_K\langle\xi\rangle^{mr},
\qquad
x\in K,\quad \|\xi\|\geq R_K;
\]
\item $a_m(x,\xi)$ is invertible for every $x\in U$ and every $\xi\neq0$;
\item there exists $b_0\in S^{-m}_{1,0}(U;\operatorname{End}(\mathbb C^r))$ such that
\[
a b_0-I_r\in S^{-1}_{1,0}(U;\operatorname{End}(\mathbb C^r)),
\qquad
b_0 a-I_r\in S^{-1}_{1,0}(U;\operatorname{End}(\mathbb C^r)).
\]
\end{enumerate}
Moreover, in the high-frequency region where $a$ is invertible, its inverse satisfies $a^{-1}\in S^{-m}_{1,0}(U;\operatorname{End}(\mathbb C^r))$.
\end{proposition}

\begin{proof}
Throughout this proof, $U,r,m$ remain fixed, and we write
\[
S^q:=S^q_{1,0}\bigl(U;\operatorname{End}(\mathbb C^r)\bigr).
\]
The identity \[
s_{\min}(a(x,\xi))=\|a(x,\xi)^{-1}\|^{-1}
\] proves equivalence of (a) and (b). If (b) holds, all singular values are bounded below by $c_K\langle\xi\rangle^m$, so their product gives (c). Conversely, the estimates for $S^m$ imply
\[
\|\operatorname{adj}a(x,\xi)\|
\leq C_K\langle\xi\rangle^{m(r-1)}.
\]
By the formula $a^{-1}=(\det a)^{-1}\operatorname{adj}a$, condition (c) yields the bound in (a).

Suppose (d) holds. On the compact set $K\times\mathbb S^{n-1}$, the function $s_{\min}(a_m)$ has a positive minimum. By definition of the classical expansion, there exists $C_K>0$ such that, uniformly for $x\in K$, $\omega\in\mathbb S^{n-1}$, and $t\geq1$,
\[
\bigl\|t^{-m}a(x,t\omega)-a_m(x,\omega)\bigr\|
\leq C_Kt^{-1}.
\]
This inequality gives (b). Conversely, if $a_m(x_0,\omega_0)v=0$ for some $v\neq0$, the same convergence contradicts the lower bound in (b) applied to $(x_0,t\omega_0)$.

If (d) holds, choose a radial cutoff $\vartheta$ vanishing near the origin and equal to one at large frequencies, and set
\[
b_0(x,\xi)=\vartheta(\xi)a_m(x,\xi)^{-1}.
\]
Homogeneity of degree $-m$ and local compactness of the sphere show that $b_0\in S^{-m}$. Since $a-\vartheta a_m\in S^{m-1}$ at high frequencies, both products in (e) are $I_r$ modulo $S^{-1}$. If (e) holds, then $ab_0=I_r+r$ with $r\in S^{-1}$. For large frequencies, $\|r\|<\displaystyle\frac{1}{2}$; thus $a$ is invertible and
\[
\|a^{-1}\|
\leq \|(I_r+r)^{-1}\|\,\|b_0\|
\leq C_K\langle\xi\rangle^{-m}.
\]

Finally, in the region where $a$ is invertible,
\[
D_{x_j}(a^{-1})=-a^{-1}(D_{x_j}a)a^{-1},
\qquad
D_{\xi_j}(a^{-1})=-a^{-1}(D_{\xi_j}a)a^{-1}.
\]
Write $\gamma=(\alpha,\beta)\in\mathbb N_0^{2n}$ and $D^\gamma=D_x^\alpha D_\xi^\beta$. Applying $D^\gamma$ to $aa^{-1}=I_r$, separating the summand in which no derivative hits $a$, and multiplying on the left by $a^{-1}$ gives the exact recurrence
\begin{equation}
\label{eq:recurrencia-derivadas-inverso-simbolo}
D^\gamma(a^{-1})
=
-a^{-1}
\sum_{0<\delta\leq\gamma}
\binom{\gamma}{\delta}
(D^\delta a)D^{\gamma-\delta}(a^{-1}),
\qquad |\gamma|\geq1.
\end{equation}
The base case $|\gamma|=0$ is the bound $\|a^{-1}\|\leq C_K\langle\xi\rangle^{-m}$ already proved. Suppose that, for every $\gamma'=(\alpha',\beta')$ with $|\gamma'|<q$,
\[
\|D^{\gamma'}a^{-1}(x,\xi)\|
\leq C_{K,\gamma'}\langle\xi\rangle^{-m-|\beta'|}.
\]
Let $|\gamma|=q$ and write $\delta=(\alpha_1,\beta_1)$. Since $0<\delta\leq\gamma$, the factor $D^{\gamma-\delta}(a^{-1})$ is covered by the induction hypothesis. Each summand in \eqref{eq:recurrencia-derivadas-inverso-simbolo} is then bounded by
\[
C_{K,\gamma,\delta}
\langle\xi\rangle^{-m}
\langle\xi\rangle^{m-|\beta_1|}
\langle\xi\rangle^{-m-(|\beta|-|\beta_1|)}
=
C_{K,\gamma,\delta}
\langle\xi\rangle^{-m-|\beta|}.
\]
The sum contains finitely many multi-indices $\delta$; thus
\[
\|D_x^\alpha D_\xi^\beta a^{-1}(x,\xi)\|
\leq C_{K,\alpha,\beta}
\langle\xi\rangle^{-m-|\beta|},
\]
which completes the induction and proves the remaining assertion.
\end{proof}

\begin{definition}[Elliptic set and characteristic set]
\label{def:conjunto-eliptico-caracteristico-pseudo}
Let $M$ be a finite-dimensional smooth manifold without boundary, and let $\mathbf{E},\mathbf{F}\to M$ be smooth complex vector bundles of the same finite rank. Let $A\in\Psi^m_{\mathrm{cl,loc}}(M;\mathbf{E},\mathbf{F})$. Its \textbf{elliptic set} is
\[
\operatorname{Ell}(A)
:=
\bigl\{(x,\xi)\in T^*M\setminus0_M
\mid \boldsymbol{\sigma}_m^\Psi(A)(x,\xi)\text{ is an isomorphism}\bigr\},
\]
and its \textbf{characteristic set} is
\[
\operatorname{Char}(A)
:=
(T^*M\setminus0_M)\setminus\operatorname{Ell}(A).
\]
The operator is \textbf{elliptic} if $\operatorname{Char}(A)=\varnothing$.
\end{definition}

This ellipticity is local in the base variable. On a noncompact manifold, it does not imply uniform ellipticity in the sense of Definition~\ref{def:simbolo-eliptico-euclidiano}. Parameter ellipticity, which requires uniform inversion of a family $\lambda I-A$ in a sector, will be defined only when constructing complex powers.

Homogeneity of the principal symbol shows that $\operatorname{Ell}(A)$ is open and conic, whereas $\operatorname{Char}(A)$ is closed and conic. In a chart, membership in $\operatorname{Ell}(A)$ is equivalent to validity of the lower bound in Proposition~\ref{prop:equivalencias-elipticidad-pseudo} on a neighborhood in the base and a conic neighborhood of the covariable. For a differential operator $P$, identity \eqref{eq:conversion-simbolos-diferencial-pseudo} shows that this definition agrees with differential ellipticity despite the factor $i^m$.

\subsection{Construction of a two-sided parametrix}

The preceding equivalences reduce ellipticity to uniform inversion of the principal component at high frequencies. Starting from this inverse, right and left corrections are constructed recursively. We first write the recurrence in coordinates, where the order of each summand can be tracked explicitly. Let $\displaystyle a\sim\displaystyle\sum_{j=0}^{\infty}a_{m-j}$ be an elliptic classical symbol. First seek an expansion \[
b^{\mathrm d}\sim\sum_{k=0}^{\infty}b^{\mathrm d}_{-m-k}
\] such that $a\# b^{\mathrm d}\sim I_r$; the letter $\mathrm d$ recalls that $\operatorname{Op}(b^{\mathrm d})$ lies to the right of $\operatorname{Op}(a)$. The degree-zero term imposes
\[
b^{\mathrm d}_{-m}=a_m^{-1}.
\]
For $N\geq1$, the component of degree $-N$ in $a\# b^{\mathrm d}$ is
\[
\sum_{\substack{j,k\geq0,\ \alpha\in\mathbb N_0^n\\
                  j+k+|\alpha|=N}}
\frac{i^{-|\alpha|}}{\alpha!}
(D_\xi^\alpha a_{m-j})
(D_x^\alpha b^{\mathrm d}_{-m-k}).
\]
Separating the unique summand containing the unknown $b^{\mathrm d}_{-m-N}$ gives the explicit recurrence
\begin{equation}
\label{eq:recurrencia-parametrica-derecha}
b^{\mathrm d}_{-m-N}
=
-a_m^{-1}
\sum_{\substack{j,k\geq0,\ \alpha\in\mathbb N_0^n\\
                  j+k+|\alpha|=N,\ k<N}}
\frac{i^{-|\alpha|}}{\alpha!}
(D_\xi^\alpha a_{m-j})
(D_x^\alpha b^{\mathrm d}_{-m-k}).
\end{equation}
To construct a left parametrix, write $b^{\mathrm i}\# a\sim I_r$. Now $b^{\mathrm i}_{-m}=a_m^{-1}$ and
\begin{equation}
\label{eq:recurrencia-parametrica-izquierda}
b^{\mathrm i}_{-m-N}
=
-\left(
\sum_{\substack{j,k\geq0,\ \alpha\in\mathbb N_0^n\\
                  j+k+|\alpha|=N,\ k<N}}
\frac{i^{-|\alpha|}}{\alpha!}
(D_\xi^\alpha b^{\mathrm i}_{-m-k})
(D_x^\alpha a_{m-j})
\right)a_m^{-1}.
\end{equation}
In particular, the first corrections are
\[
\begin{split}
b^{\mathrm d}_{-m-1}
={}&-a_m^{-1}\left(
a_{m-1}a_m^{-1}
+i^{-1}\sum_{\ell=1}^n
(D_{\xi_\ell}a_m)(D_{x_\ell}a_m^{-1})
\right),\\
b^{\mathrm i}_{-m-1}
={}&-\left(
a_m^{-1}a_{m-1}
+i^{-1}\sum_{\ell=1}^n
(D_{\xi_\ell}a_m^{-1})(D_{x_\ell}a_m)
\right)a_m^{-1}.
\end{split}
\]
The order of factors matters for matrix-valued symbols. The factors $i^{-|\alpha|}$ likewise cannot be omitted: they come exactly from the convention $D=\partial$ fixed at the beginning of the chapter.

\begin{theorem}[Existence and equivalence of parametrices]
\label{teo:parametrica-eliptica-pseudo}
Let $M$ be a finite-dimensional smooth manifold without boundary, and let $\mathbf{E},\mathbf{F}\to M$ be smooth complex vector bundles of the same finite rank. Let $A\in\Psi^m_{\mathrm{cl}}(M;\mathbf{E},\mathbf{F})$. The following conditions are equivalent:
\begin{enumerate}[label=(\alph*)]
\item $A$ is elliptic;
\item there exists $B\in\Psi^{-m}_{\mathrm{cl}}(M;\mathbf{F},\mathbf{E})$ such that
\[
AB=I_{\mathbf{F}}-R_{\mathbf{F}},
\qquad
BA=I_{\mathbf{E}}-R_{\mathbf{E}},
\qquad
R_{\mathbf{E}}\in\Psi^{-\infty}(M;\mathbf{E},\mathbf{E}),
\qquad
R_{\mathbf{F}}\in\Psi^{-\infty}(M;\mathbf{F},\mathbf{F});
\]
\item the principal symbol of $A$ has a homogeneous inverse of degree $-m$ on $T^*M\setminus0_M$.
\end{enumerate}
The operator $B$ can be chosen properly supported and is a two-sided parametrix of $A$. Two two-sided parametrices differ by an operator in $\Psi^{-\infty}(M;\mathbf{F},\mathbf{E})$.
\end{theorem}

\begin{proof}
On an open set $U\subseteq\mathbb R^n$ and in frames of common rank $r$, formulas \eqref{eq:recurrencia-parametrica-derecha} and \eqref{eq:recurrencia-parametrica-izquierda} define homogeneous terms of degrees $-m-k$. Let us verify this inductively. For $k=0$, $a_m^{-1}(x,t\xi)=t^{-m}a_m(x,\xi)^{-1}$. Suppose $b^{\mathrm d}_{-m-k}$ is homogeneous of degree $-m-k$ for $k<N$. In each summand of \eqref{eq:recurrencia-parametrica-derecha}, $D_\xi^\alpha a_{m-j}$ has degree $m-j-|\alpha|$ and $D_x^\alpha b^{\mathrm d}_{-m-k}$ retains degree $-m-k$; since $j+k+|\alpha|=N$, their product has degree $-N$. Multiplication by $a_m^{-1}$ lowers the degree by $m$ and gives exactly $-m-N$. The sum is finite because $j,k,|\alpha|\leq N$. The same occurs in the left recurrence: the sum in parentheses has degree $-N$, and the factor $a_m^{-1}$ on the right gives degree $-m-N$. This proves both inductive steps simultaneously. Lemma~\ref{lem:suma-asintotica-simbolos} produces symbols $b^{\mathrm d},b^{\mathrm i}\in
S^{-m}_{\mathrm{cl}}(U;\operatorname{End}(\mathbb C^r))$ with these expansions. Theorem~\ref{teo:composicion-simbolica-pseudo} and the recurrences give
\[
a\# b^{\mathrm d}-I_r
\in S^{-\infty}(U\times\mathbb R^n;
\operatorname{End}(\mathbb C^r)),
\qquad
b^{\mathrm i}\# a-I_r
\in S^{-\infty}(U\times\mathbb R^n;
\operatorname{End}(\mathbb C^r)).
\]
Quantizing and arranging proper support yields a right parametrix $B_{\mathrm d}$ and a left parametrix $B_{\mathrm i}$. If $AB_{\mathrm d}=I_{\mathbf{F}}-R_{\mathrm d}$ and $B_{\mathrm i}A=I_{\mathbf{E}}-R_{\mathrm i}$, then
\[
B_{\mathrm i}-B_{\mathrm d}
=B_{\mathrm i}(I_{\mathbf{F}}-AB_{\mathrm d})
 -(I_{\mathbf{E}}-B_{\mathrm i}A)B_{\mathrm d}
=B_{\mathrm i}R_{\mathrm d}-R_{\mathrm i}B_{\mathrm d}.
\]
Corollary~\ref{cor:ideal-regularizantes-pseudo} shows that this difference is smoothing. Thus either parametrix is two-sided modulo $\Psi^{-\infty}(M;\mathbf{E},\mathbf{E})$ in the left defect and $\Psi^{-\infty}(M;\mathbf{F},\mathbf{F})$ in the right defect.

Globally, (c) supplies the homogeneous section $\mathbf{a}_m^{-1}:=\boldsymbol{\sigma}_m^\Psi(A)^{-1}$. Surjectivity of the exact sequence in Theorem~\ref{teo:sucesion-exacta-simbolo-principal-pseudo} allows it to be quantized as an operator $B_0$ of order $-m$; then $AB_0-I_{\mathbf{F}}$ and $B_0A-I_{\mathbf{E}}$ have order $-1$. Set $B_0^{\mathrm d}=B_0^{\mathrm i}=B_0$.

The right induction statement is
\[
B_{<N}^{\mathrm d}:=\sum_{k=0}^{N-1}B_k^{\mathrm d},
\qquad
B_k^{\mathrm d}\in\Psi^{-m-k}_{\mathrm{cl}}(M;\mathbf{F},\mathbf{E}),
\qquad
AB_{<N}^{\mathrm d}-I_{\mathbf{F}}\in\Psi^{-N}_{\mathrm{cl}}(M;\mathbf{F},\mathbf{F}).
\]
The base case $N=1$ is the choice of $B_0$. Assuming the assertion for $N$, define the global homogeneous section
\[
\mathbf{e}_N^{\mathrm d}
:=
\boldsymbol{\sigma}_{-N}^\Psi\bigl(AB_{<N}^{\mathrm d}-I_{\mathbf{F}}\bigr)
\in
\mathcal S^{-N}_{\mathrm{hom}}
\bigl(T^*M;\operatorname{End}(\mathbf{F})\bigr)
\]
and
\[
\mathbf{c}_{-m-N}^{\mathrm d}:=-\mathbf{a}_m^{-1}\mathbf{e}_N^{\mathrm d}
\in
\mathcal S^{-m-N}_{\mathrm{hom}}
\bigl(T^*M;\operatorname{Hom}(\mathbf{F},\mathbf{E})\bigr).
\]
These are global sections: $\mathbf{e}_N^{\mathrm d}$ is the principal symbol of a global operator, and the bundle composition $\mathbf{a}_m^{-1}\circ \mathbf{e}_N^{\mathrm d}$ transforms $\pi^*\mathbf{F}$ into $\pi^*\mathbf{E}$. Surjectivity of the exact sequence provides $B_N^{\mathrm d}\in\Psi^{-m-N}_{\mathrm{cl}}(M;\mathbf{F},\mathbf{E})$ with principal symbol $\mathbf{c}_{-m-N}^{\mathrm d}$. The product rule for principal symbols gives
\[
\begin{split}
\boldsymbol{\sigma}_{-N}^\Psi
\bigl(A(B_{<N}^{\mathrm d}+B_N^{\mathrm d})-I_{\mathbf{F}}\bigr)
&=\mathbf{e}_N^{\mathrm d}+\boldsymbol{\sigma}_m^\Psi(A)\mathbf{c}_{-m-N}^{\mathrm d}\\
&=\mathbf{e}_N^{\mathrm d}-\boldsymbol{\sigma}_m^\Psi(A)\mathbf{a}_m^{-1}\mathbf{e}_N^{\mathrm d}=0.
\end{split}
\]
By exactness, the new error belongs to $\Psi^{-N-1}_{\mathrm{cl}}(M;\mathbf{F},\mathbf{F})$; this completes the right step.

For the left side, suppose $B_{<N}^{\mathrm i}A-I_{\mathbf{E}}\in\Psi^{-N}_{\mathrm{cl}}(M;\mathbf{E},\mathbf{E})$ and set
\[
\mathbf{e}_N^{\mathrm i}
:=
\boldsymbol{\sigma}_{-N}^\Psi\bigl(B_{<N}^{\mathrm i}A-I_{\mathbf{E}}\bigr)
\in
\mathcal S^{-N}_{\mathrm{hom}}
\bigl(T^*M;\operatorname{End}(\mathbf{E})\bigr),
\qquad
\mathbf{c}_{-m-N}^{\mathrm i}:=-\mathbf{e}_N^{\mathrm i}\mathbf{a}_m^{-1}.
\]
Now $\mathbf{c}_{-m-N}^{\mathrm i}$ is a global section of $\operatorname{Hom}(\pi^*\mathbf{F},\pi^*\mathbf{E})$ and can be quantized as $B_N^{\mathrm i}\in\Psi^{-m-N}_{\mathrm{cl}}(M;\mathbf{F},\mathbf{E})$. In this case,
\[
\boldsymbol{\sigma}_{-N}^\Psi
\bigl((B_{<N}^{\mathrm i}+B_N^{\mathrm i})A-I_{\mathbf{E}}\bigr)
=\mathbf{e}_N^{\mathrm i}+\mathbf{c}_{-m-N}^{\mathrm i}\boldsymbol{\sigma}_m^\Psi(A)
=0,
\]
so the left error also drops from $-N$ to $-N-1$. In a chart, the components of the two choices of $\mathbf{c}_{-m-N}$ are exactly the right-hand sides of \eqref{eq:recurrencia-parametrica-derecha} and \eqref{eq:recurrencia-parametrica-izquierda}; the preceding global formulas prove, rather than presuppose, that these expressions glue.

Lemma~\ref{lem:suma-asintotica-simbolos}, applied over a locally finite cover and followed by adjustment of proper support, sums the two families into operators $B_{\mathrm d}$ and $B_{\mathrm i}$. The inductive steps show that $AB_{\mathrm d}-I_{\mathbf{F}}$ belongs to $\Psi^{-N}(M;\mathbf{F},\mathbf{F})$ and $B_{\mathrm i}A-I_{\mathbf{E}}$ to $\Psi^{-N}(M;\mathbf{E},\mathbf{E})$ for each $N$, so both are smoothing. The algebraic identity above shows that $B_{\mathrm i}-B_{\mathrm d}\in\Psi^{-\infty}(M;\mathbf{F},\mathbf{E})$ and that either operator is two-sided. This reproduces the local recurrences globally and proves (c)$\Rightarrow$(b). The implication (b)$\Rightarrow$(c) follows by taking principal symbols in $AB=I_{\mathbf{F}}-R_{\mathbf{F}}$ and $BA=I_{\mathbf{E}}-R_{\mathbf{E}}$. Equivalence of (a) and (c) is Proposition~\ref{prop:equivalencias-elipticidad-pseudo} applied in charts. If $B$ and $\widetilde B$ are two-sided, the same identity above, with $B_{\mathrm i}=B$ and $B_{\mathrm d}=\widetilde B$, shows that $B-\widetilde B$ is smoothing.
\end{proof}

\begin{remark}[A parametrix need not be an inverse]
\label{obs:parametrica-no-inversa-pseudo}
The theorem's identities hold modulo smooth kernels, not exactly. A parametrix exists because of high-frequency behavior; an inverse also requires injectivity and surjectivity in the space where $A$ is realized. For example, on a closed connected manifold, the Laplacian is elliptic and has a parametrix, but constant functions belong to its kernel. It is therefore not invertible on all of $C^\infty(M)$. If a realization of $A$ has a true inverse, that inverse agrees with any parametrix modulo the correction eliminating the smoothing remainders; this additional assertion is not part of symbolic ellipticity.
\end{remark}

\begin{corollary}[Qualitative elliptic regularity]
\label{cor:regularidad-eliptica-cualitativa-pseudo}
Let $M$ be a finite-dimensional smooth manifold without boundary, and let $\mathbf{E},\mathbf{F}\to M$ be smooth complex vector bundles of the same finite rank. Let $A\in\Psi^m_{\mathrm{cl}}(M;\mathbf{E},\mathbf{F})$ be elliptic and properly supported. If $\mathbf{u}\in\mathcal D'(M,\mathbf{E})$ and $A\mathbf{u}$ is smooth on an open set $V\subseteq M$, then $\mathbf{u}$ is smooth on $V$. In particular,
\[
\operatorname{sing\,supp}(A\mathbf{u})
=
\operatorname{sing\,supp}(\mathbf{u}).
\]
\end{corollary}

\begin{proof}
Let $B$ be a two-sided parametrix. From $BA=I_{\mathbf{E}}-R_{\mathbf{E}}$ we obtain $\mathbf{u}=B(A\mathbf{u})+R_{\mathbf{E}}\mathbf{u}$. The first summand is smooth on $V$ by pseudolocality of $B$, and the second is smooth on all of $M$. This proves the first assertion. The reverse inclusion for singular supports is Theorem~\ref{teo:pseudolocalidad-pseudo}.
\end{proof}

\subsection{The wavefront set on a Euclidean open set}

Qualitative elliptic regularity detects the points where a distribution is smooth. To also record cotangent directions in which smoothness fails, we must localize in the base and examine Fourier decay within cones.

Write
\[
\mathcal E'(\mathbb R^n)
:=
\{v\in\mathcal D'(\mathbb R^n)\mid \supp v\text{ is compact}\}.
\]
A compactly supported distribution $v$ is tempered, and its Fourier transform is the smooth function
\[
\widehat v(\xi)
=(2\pi)^{-\frac{n}{2}}\langle v(x),e^{-ix\cdot\xi}\rangle.
\]
Continuity of $v$ on functions supported in a fixed compact set shows that $\widehat v$ and all its derivatives have at most polynomial growth. The wavefront set records the directions where this growth improves to rapid decay after localization in the spatial variable.

We say that two open cones $\Gamma',\Gamma\subseteq\mathbb R^n\setminus
\{0\}$ satisfy $\Gamma'\Subset_{\mathrm c}\Gamma$ if
\[
\overline{\Gamma'\cap\mathbb S^{n-1}}
\subseteq\Gamma\cap\mathbb S^{n-1}.
\]

\begin{lemma}[Stability under shrinking the cone]
\label{lem:estabilidad-conica-fourier}
Let $v\in\mathcal E'(\mathbb R^n)$ and suppose that, for every $N$, there exists $C_N$ such that
\[
|\widehat v(\boldsymbol{\xi})|\leq C_N\langle\boldsymbol{\xi}\rangle^{-N},
\qquad \boldsymbol{\xi}\in\Gamma,
\]
where $\Gamma$ is an open cone. If $\psi\in C_c^\infty(\mathbb R^n)$ and $\Gamma'\Subset_{\mathrm c}\Gamma$, then $\widehat{\psi v}$ decays rapidly in $\Gamma'$.
\end{lemma}

\begin{proof}
With the unitary normalization,
\[
\widehat{\psi v}(\boldsymbol{\xi})
=(2\pi)^{-\frac{n}{2}}
\int_{\mathbb R^n}\widehat\psi(\boldsymbol{\xi}-\boldsymbol{\eta})\widehat v(\boldsymbol{\eta})\,d\boldsymbol{\eta}.
\]
The function $\widehat v$ has polynomial growth, and $\widehat\psi$ is a Schwartz function. For $\boldsymbol{\xi}\in\Gamma'$, split the integral into $\boldsymbol{\eta}\in\Gamma$ and $\boldsymbol{\eta}\notin\Gamma$. In the first region, use the arbitrary decay of $\widehat v$ and Peetre's inequality in Lemma~\ref{lem:desigualdad-peetre-peso-japones}. In the second, angular separation of the cones gives a constant $c>0$ such that \[
|\boldsymbol{\xi}-\boldsymbol{\eta}|\geq c(|\boldsymbol{\xi}|+|\boldsymbol{\eta}|)
\] if $\boldsymbol{\xi}\in\Gamma'$ and $\boldsymbol{\eta}\notin\Gamma$. The decay of $\widehat\psi$ then offsets any power of $\boldsymbol{\xi}$ as well as the polynomial growth of $\widehat v$. Both integrals are bounded by $C_N'\langle\boldsymbol{\xi}\rangle^{-N}$ for every $N$.
\end{proof}

\begin{definition}[Euclidean wavefront set]
\label{def:frente-onda-euclidiano}
\index{wavefront set}
\glsadd{frente-onda}
Let $U\subseteq\mathbb R^n$ be open and let $u\in\mathcal D'(U)$. A point $(x_0,\boldsymbol{\xi}_0)\in U\times(\mathbb R^n\setminus\{0\})$ is a \textbf{regular direction} of $u$ if there exist $\varphi\in C_c^\infty(U)$ equal to one on a neighborhood of $x_0$ and an open cone $\Gamma$ containing $\boldsymbol{\xi}_0$ such that, for every $N\in\mathbb
N_0$,
\[
\sup_{\boldsymbol{\xi}\in\Gamma}
\langle\boldsymbol{\xi}\rangle^N|\widehat{\varphi u}(\boldsymbol{\xi})|<\infty.
\]
The \textbf{wavefront set} of $u$ is the complement of the regular directions:
\[
\operatorname{WF}(u)
\subseteq U\times(\mathbb R^n\setminus\{0\}).
\]
\end{definition}

Lemma~\ref{lem:estabilidad-conica-fourier} shows that the definition is independent of the cutoff: if $\psi$ is supported in the region where $\varphi=1$, then $\psi u=\psi(\varphi u)$, and it suffices to shrink the cone.

\begin{proposition}[Elementary properties of the wavefront set]
\label{prop:propiedades-elementales-wf}
Let $U\subseteq\mathbb R^n$ be open. For $u,v\in\mathcal D'(U)$ and $f\in C^\infty(U)$, the following hold:
\begin{enumerate}[label=(\roman*)]
\item $\operatorname{WF}(u)$ is closed in $T^*U\setminus0_U$ and conic in the covariable;
\item \[
\operatorname{WF}(u+v)
\subseteq\operatorname{WF}(u)\cup\operatorname{WF}(v),
\qquad
\operatorname{WF}(fu)\subseteq\operatorname{WF}(u);
\]
\item the projection $\pi\colon T^*U\longrightarrow U$ satisfies
\[
\pi\bigl(\operatorname{WF}(u)\bigr)
=\operatorname{sing\,supp}(u).
\]
\end{enumerate}
\end{proposition}

\begin{proof}
If a cutoff and a cone establish regularity at $(x_0,\xi_0)$, the same function works for base points in a neighborhood of $x_0$, and the same cone works for directions near $\xi_0$. Thus the set of regular directions is open. Invariance under positive dilations is immediate and proves (i).

For a sum, choose a common cutoff and intersect two cones of regularity. For $fu$, take a cutoff $\psi$ supported where the original cutoff equals one, and apply Lemma~\ref{lem:estabilidad-conica-fourier} to the smooth multiplier $\psi f$. This proves (ii).

If $u$ is smooth near $x_0$, a localization $\varphi u$ belongs to $C_c^\infty$ and its transform is Schwartz in every direction. Thus $x_0$ does not belong to the projection of the wavefront set. Conversely, suppose no nonzero covariable over $x_0$ belongs to the wavefront set. Compactness of $\mathbb S^{n-1}$ allows finitely many cones of regularity to be chosen covering the sphere. Choose a single function $\varphi$ equal to one near $x_0$ and supported in the intersection of the regions where the corresponding cutoffs equal one. The preceding lemma then shows that $\widehat{\varphi u}$ decays faster than any power in every direction. Thus $\xi^\alpha\widehat{\varphi u}\in L^1(\mathbb R^n)$ for every $\alpha$. Fourier inversion in the distributional sense implies $\varphi u\in C^\infty$. Hence $u$ is smooth near $x_0$, proving (iii).
\end{proof}

\subsection{Invariance and the global definition}

The Euclidean definition is expressed in frequencies, but the direction it records belongs to the cotangent bundle. We must check that a diffeomorphism transports these directions by its cotangent map and preserves rapid decay in the corresponding cones. The precise assertion is as follows.

\begin{proposition}[Invariance under diffeomorphisms]
\label{prop:invariancia-difeomorfismos-wf}
Let $\kappa\colon U\longrightarrow V$ be a diffeomorphism between Euclidean open sets, and let $u\in\mathcal D'(U)$. For the transformed distribution $\kappa_*u\in\mathcal D'(V)$,
\[
\operatorname{WF}(\kappa_*u)
=
\bigl\{\bigl(\kappa(x),D\kappa(x)^{-T}\xi\bigr)
\mid (x,\xi)\in\operatorname{WF}(u)\bigr\}.
\]
\end{proposition}

\begin{proof}
Localize $u$ on a compact set where $\kappa$ is the chart under consideration. By Fourier inversion of the localized distribution, the Fourier transform of a localization of $\kappa_*u$ can be written as an integral in $(x,\xi)$ with phase
\[
x\cdot\xi-\kappa(x)\cdot\eta.
\]
Its stationary points in $x$ satisfy $\xi=D\kappa(x)^T\eta$. Suppose $(x_0,\xi_0)$ is regular, and set $\eta_0=D\kappa(x_0)^{-T}\xi_0$. Shrink the neighborhood of $x_0$ and the cones of $\xi_0,\eta_0$ so that every solution of the stationary relation remains within the cone where the localized transform of $u$ decays rapidly. Use this decay in that region. In the complement, the gradient $\xi-D\kappa(x)^T\eta$ is bounded away from zero; repeated integration by parts with
\[
L=
\frac{1-i(\xi-D\kappa(x)^T\eta)\cdot D_x}{1+\|\xi-D\kappa(x)^T\eta\|^2}
\]
produces any negative power of $\langle\eta\rangle$. Jacobians and cutoffs are smooth compactly supported amplitudes and do not change the estimates. Thus $(\kappa(x_0),\eta_0)$ is regular for $\kappa_*u$. Applying the same argument to $\kappa^{-1}$ gives the reverse inclusion and proves equality.
\end{proof}

\begin{definition}[Global wavefront set]
\label{def:frente-onda-global}
Let $M$ be a finite-dimensional smooth manifold without boundary, let $\mathbf{E}\to M$ be a smooth complex vector bundle of finite rank, and let $\mathbf{u}\in\mathcal D'(M,\mathbf{E})$. In a chart $\kappa\colon U\to\mathbb R^n$ and a local frame, write $\mathbf{u}=(u^1,\dots,u^r)$. Declare that $(x,\xi)\in T^*M\setminus0_M$ does not belong to $\operatorname{WF}(\mathbf{u})$ if
\[
\bigl(\kappa(x),D\kappa(x)^{-T}\xi\bigr)
\notin\bigcup_{a=1}^r\operatorname{WF}(u^a).
\]
\end{definition}

Proposition~\ref{prop:invariancia-difeomorfismos-wf} proves independence of the chart. A frame change multiplies the components by a smooth invertible matrix; Proposition~\ref{prop:propiedades-elementales-wf}(ii), applied in both directions, proves independence of the frame. Thus $\operatorname{WF}(\mathbf{u})$ is a closed conic subset of $T^*M\setminus0_M$, and
\[
\pi\bigl(\operatorname{WF}(\mathbf{u})\bigr)
=\operatorname{sing\,supp}(\mathbf{u}).
\]
The definition is also independent of the choice of a positive density for representing distributions in coordinates: two choices differ by a nonvanishing smooth function. In the Riemannian case, as elsewhere in Riemannian analysis, we use the measure $d\lambda_{\mathbf{g}}$. In this precise sense, the wavefront set is canonical: it is independent of the atlas, frame, metric, and auxiliary density.

\subsection{Pseudodifferential action and microlocal regularity}

Pseudolocality says that a pseudodifferential operator creates no singularities in the base. The wavefront set refines this assertion: the operator creates no new singular directions either. When the symbol is elliptic at a covector, a microlocal parametrix also recovers the datum's singular direction from that of its image.

\begin{theorem}[Pseudodifferential operators create no wavefronts]
\label{teo:wf-accion-pseudodiferencial}
Let $M$ be a finite-dimensional smooth manifold without boundary, and let $\mathbf{E},\mathbf{F}\to M$ be smooth complex vector bundles of finite rank. Let $A\in\Psi^m(M;\mathbf{E},\mathbf{F})$ be properly supported. Then, for every $\mathbf{u}\in\mathcal D'(M,\mathbf{E})$,
\[
\operatorname{WF}(A\mathbf{u})
\subseteq
\operatorname{WF}(\mathbf{u}).
\]
In particular, for every multi-index $\alpha$ and every smooth function $f$,
\[
\operatorname{WF}(D^\alpha \mathbf{u})\subseteq\operatorname{WF}(\mathbf{u}),
\qquad
\operatorname{WF}(f\mathbf{u})\subseteq\operatorname{WF}(\mathbf{u})
\]
in coordinates.
\end{theorem}

\begin{proof}
The assertion is local. In coordinates, fix $(x_0,\xi_0)\notin\operatorname{WF}(\mathbf{u})$ and write $A=\operatorname{Op}(a)$ modulo a smoothing operator. Choose $\psi\in C_c^\infty$ equal to one near $x_0$ and supported where a cutoff establishing regularity of $(x_0,\xi_0)$ equals one, and decompose $\mathbf{u}=\psi \mathbf{u}+(1-\psi)\mathbf{u}$. The second summand yields a smooth function near $x_0$ because the kernel of $A$ is smooth off the diagonal. For $\mathbf{v}=\psi \mathbf{u}\in\mathcal E'$ and $\varphi\in C_c^\infty$ supported near $x_0$, the quantization formula gives
\begin{equation}
\label{eq:fourier-localizacion-Au}
\widehat{\varphi A\mathbf{v}}(\zeta)
=(2\pi)^{-\frac{n}{2}}
\int_{\mathbb R^n}
\mathcal F_{x\to\theta}
\bigl(\varphi(x)a(x,\boldsymbol{\eta})\bigr)(\zeta-\boldsymbol{\eta})
\widehat{\mathbf{v}}(\boldsymbol{\eta})\,d\boldsymbol{\eta}.
\end{equation}
Integrating by parts in $x$ gives, for every $L$,
\[
\left\|
\mathcal F_{x\to\theta}
\bigl(\varphi a(\cdot,\boldsymbol{\eta})\bigr)(\theta)
\right\|
\leq C_L\langle\theta\rangle^{-L}\langle\boldsymbol{\eta}\rangle^m.
\]
Take a cone $\Gamma$ where $\widehat{\mathbf{v}}$ decays rapidly and a cone $\Gamma'\Subset_{\mathrm c}\Gamma$ containing $\xi_0$. In \eqref{eq:fourier-localizacion-Au}, separate $\boldsymbol{\eta}\in\Gamma$ and $\boldsymbol{\eta}\notin\Gamma$. In the first part, decay of $\widehat{\mathbf{v}}$ absorbs the order $m$; in the second, use $|\zeta-\boldsymbol{\eta}|\geq c(|\zeta|+|\boldsymbol{\eta}|)$ for $\zeta\in\Gamma'$. Since $L$ is arbitrary and $\widehat{\mathbf{v}}$ has polynomial growth, the right-hand side decays rapidly in $\Gamma'$. This proves the inclusion. The final two assertions are special cases, since derivatives and smooth multiplication are pseudodifferential operators.
\end{proof}

To invert only near a covector, it is useful to fix some terminology. Let $M$ be a finite-dimensional smooth manifold without boundary, and let $\mathbf{E},\mathbf{F}\longrightarrow M$ be smooth complex vector bundles of finite rank. We say that a properly supported $R\in\Psi^s(M;\mathbf{E},\mathbf{F})$ is \textbf{microlocally smoothing} at $q=(x_0,\xi_0)\in T^*M\setminus0_M$ if, in some chart and local frames of $\mathbf{E}$ and $\mathbf{F}$, the symbol of a localization of $R$ satisfies, on a product $V\times\Gamma$ containing $q$,
\[
\|D_x^\alpha D_\xi^\beta r(x,\xi)\|
\leq C_{N,\alpha,\beta}\langle\xi\rangle^{-N}
\]
for every $N$ and large frequencies. The proof of the preceding theorem, using these bounds in the corresponding conic part and cone separation in its complement, gives
\begin{equation}
\label{eq:microlocalmente-regularizante-wf}
q\notin\operatorname{WF}(R\mathbf{u})
\qquad\text{for every }\mathbf{u}\in\mathcal D'(M,\mathbf{E}).
\end{equation}
Let us check that this condition is independent of the chart. Let $\kappa\colon U\to\widetilde U$ be a coordinate change, write \[
\mathcal B_\kappa(x,y)
:=
\int_0^1D\kappa\bigl(y+t(x-y)\bigr)\,dt
\], and suppose $r$ satisfies the preceding bounds near $(\kappa(x_0),\eta_0)$, where $\eta_0=D\kappa(x_0)^{-T}\xi_0$. After inserting a cutoff $\rho(x,y)$ equal to one near the diagonal, the formula in Theorem~\ref{teo:cambio-coordenadas-pseudo} gives the transformed amplitude
\begin{equation}
\label{eq:amplitud-cambio-carta-microlocal}
\widetilde b(x,y,\xi)
=
\rho(x,y)
r\bigl(\kappa(x),\mathcal B_\kappa(x,y)^{-T}\xi\bigr)
\frac{|\det D\kappa(y)|}{|\det\mathcal B_\kappa(x,y)|}.
\end{equation}
Frame changes merely add smooth matrices on the left and right of this expression.

By continuity of $\mathcal B_\kappa^{-T}$ on unit directions, one can choose cones $\Gamma_2\Subset_{\mathrm c}\Gamma_1$ around $\boldsymbol{\xi}_0$ and a neighborhood $U_0$ of $x_0$ such that
\[
\mathcal B_\kappa(x,y)^{-T}\boldsymbol{\xi}\in\Gamma,
\qquad
c\langle\boldsymbol{\xi}\rangle
\leq
\langle\mathcal B_\kappa(x,y)^{-T}\boldsymbol{\xi}\rangle
\leq C\langle\boldsymbol{\xi}\rangle
\]
for $x,y\in U_0$, $\boldsymbol{\xi}\in\Gamma_1$, and large frequencies. Fix an angular cutoff $\vartheta$ supported in $\Gamma_1$ and equal to one on $\Gamma_2$. The chain rule shows that each derivative of $\vartheta(\boldsymbol{\xi})\widetilde b(x,y,\boldsymbol{\xi})$ is a finite sum of terms
\[
c_{\mu,\nu,\lambda}(x,y)\,\xi^\mu
(D_z^\nu D_\eta^\lambda r)
\bigl(\kappa(x),\mathcal B_\kappa(x,y)^{-T}\boldsymbol{\xi}\bigr),
\]
and terms of the same form with a derivative of $\vartheta$. The coefficients $c_{\mu,\nu,\lambda}$ and all their derivatives are bounded after shrinking $U_0$. Given $L$, apply the bound on $r$ with exponent $L+|\mu|$; the preceding comparability gives, for all $\alpha,\gamma,\beta$,
\begin{equation}
\label{eq:cota-amplitud-cambio-carta-microlocal}
\|D_x^\alpha D_y^\gamma D_\xi^\beta
(\vartheta\widetilde b)(x,y,\boldsymbol{\xi})\|
\leq
C_{L,\alpha,\gamma,\beta}\langle\boldsymbol{\xi}\rangle^{-L}
\end{equation}
on $U_0\times U_0\times\Gamma_2$.

When reducing \eqref{eq:amplitud-cambio-carta-microlocal}, each term \[
\frac{i^{-|\delta|}}{\delta!}
D_\xi^\delta D_y^\delta\widetilde b(x,y,\xi)\big|_{y=x}
\] satisfies \eqref{eq:cota-amplitud-cambio-carta-microlocal} for any $L$. If $s$ is a global order of $r$, the remainder after summing all terms with $|\delta|<J$ satisfies symbol estimates of order $s-J$ by Theorem~\ref{teo:amplitud-a-simbolo}. For a prescribed bound $\langle\xi\rangle^{-L}$, choose $J>s+L$; the remainder then also satisfies that bound in the smaller cone. Thus the transformed symbol satisfies estimates of order $-L$ on $U_0\times\Gamma_2$ for every $L$. Applying the same calculation to $\kappa^{-1}$ proves equivalence of the two conditions. The smooth frame and density factors already considered change none of the bounds. This proves independence of the chart and auxiliary trivializations.

\begin{proposition}[Microlocal parametrix]
\label{prop:parametrica-microlocal-pseudo}
Let $M$ be a finite-dimensional smooth manifold without boundary, and let $\mathbf{E},\mathbf{F}\to M$ be smooth complex vector bundles of the same finite rank. If $q\in\operatorname{Ell}(A)$ for a properly supported $A\in\Psi^m_{\mathrm{cl}}(M;\mathbf{E},\mathbf{F})$, there exist a properly supported $B\in\Psi^{-m}_{\mathrm{cl}}(M;\mathbf{F},\mathbf{E})$ and properly supported operators $R_{\mathbf{E}},R_{\mathbf{F}}$ microlocally smoothing at $q$ such that
\[
BA=I_{\mathbf{E}}-R_{\mathbf{E}},
\qquad
AB=I_{\mathbf{F}}-R_{\mathbf{F}}
\]
microlocally at $q$. For $A\in\Psi^m_{\mathrm{cl,loc}}$, the same
conclusion holds for a properly supported representative obtained by a
cutoff near the diagonal; it differs from $A$ by a smooth kernel.
\end{proposition}

\begin{proof}
Choose a chart $\kappa\colon U\longrightarrow\Omega\subseteq\mathbb R^n$ and local frames of common rank $r$. In this proof write
\[
S^j:=S^j_{1,0}
\bigl(\Omega;\operatorname{End}(\mathbb C^r)\bigr).
\]
Choose neighborhoods $V_2\Subset V_1\Subset\Omega$ of the base point and cones $\Gamma_2\Subset_{\mathrm c}\Gamma_1$ around the covariable of $q$ so that $a_m$ is invertible on $V_1\times\Gamma_1$. Let $\chi(x,\xi)$ be a cutoff of degree zero at large frequencies, supported in $V_1\times\Gamma_1$ and equal to one on $V_2\times\Gamma_2$. Set
\[
b_{-m}^{\mathrm d}=b_{-m}^{\mathrm i}=\chi a_m^{-1}.
\]
Since all derivatives of $\chi$ vanish on the smaller product, there $a\# b_{-m}^{\mathrm d}-I_{\mathbf{F}}\in S^{-1}$ and $b_{-m}^{\mathrm i}\#a-I_{\mathbf{E}}\in S^{-1}$.

The right induction statement is that, after constructing $b_{-m-k}^{\mathrm d}$ for $0\leq k<N$,
\[
a\#\sum_{k=0}^{N-1}b_{-m-k}^{\mathrm d}-I_{\mathbf{F}}
\in S^{-N}
\qquad\text{on }V_2\times\Gamma_2.
\]
Let $e_N^{\mathrm d}$ be the error's homogeneous component of degree $-N$, computed on $V_1\times\Gamma_1$. Define
\[
b_{-m-N}^{\mathrm d}=-\chi a_m^{-1}e_N^{\mathrm d}.
\]
On $V_2\times\Gamma_2$, the new component of degree $-N$ is $e_N^{\mathrm d}+a_mb_{-m-N}^{\mathrm d}=0$; the composition formula therefore gives an error in $S^{-N-1}$. This is exactly the step in \eqref{eq:recurrencia-parametrica-derecha} and proves the right induction.

For the left side, if $e_N^{\mathrm i}$ is the component of degree $-N$ of \[
\left(\sum_{k=0}^{N-1}b_{-m-k}^{\mathrm i}\right)\#a-I_{\mathbf{E}},
\] computed on $V_1\times\Gamma_1$, take $b_{-m-N}^{\mathrm i}=-\chi e_N^{\mathrm i}a_m^{-1}$. The equality $e_N^{\mathrm i}+b_{-m-N}^{\mathrm i}a_m=0$ again reduces the error to $S^{-N-1}$ and reproduces \eqref{eq:recurrencia-parametrica-izquierda}.

Lemma~\ref{lem:suma-asintotica-simbolos} sums both families into symbols $b^{\mathrm d},b^{\mathrm i}\in S^{-m}$. For each $N$, both errors belong to $S^{-N}$ on $V_2\times\Gamma_2$, so they are microlocally smoothing at $q$. At the operator level,
\[
B_{\mathrm i}-B_{\mathrm d}
=B_{\mathrm i}(I_{\mathbf{F}}-AB_{\mathrm d})
 -(I_{\mathbf{E}}-B_{\mathrm i}A)B_{\mathrm d}.
\]
The composition expansion shows that composing a finite-order symbol with one belonging to $S^{-N}$ on the smaller product for every $N$ preserves this property, after shrinking the cones once more. Thus $B_{\mathrm i}-B_{\mathrm d}$ is microlocally smoothing, and $B:=B_{\mathrm d}$ also satisfies the left identity. Quantizing, inserting spatial cutoffs, and arranging proper support yields the two asserted identities; these changes add only smooth kernels near $q$.
\end{proof}

\begin{theorem}[Microlocal elliptic regularity]
\label{teo:regularidad-eliptica-microlocal}
Let $M$ be a finite-dimensional smooth manifold without boundary, and let $\mathbf{E},\mathbf{F}\to M$ be smooth complex vector bundles of the same finite rank. Let $A\in\Psi^m_{\mathrm{cl}}(M;\mathbf{E},\mathbf{F})$ be properly supported. For every $\mathbf{u}\in\mathcal D'(M,\mathbf{E})$,
\[
\operatorname{WF}(A\mathbf{u})
\subseteq\operatorname{WF}(\mathbf{u})
\subseteq
\operatorname{WF}(A\mathbf{u})\cup\operatorname{Char}(A).
\]
Equivalently,
\begin{equation}
\label{eq:igualdad-eliptica-wf}
\operatorname{WF}(A\mathbf{u})\cap\operatorname{Ell}(A)
=
\operatorname{WF}(\mathbf{u})\cap\operatorname{Ell}(A).
\end{equation}
If $A$ is elliptic, then $\operatorname{WF}(A\mathbf{u})=\operatorname{WF}(\mathbf{u})$.
\end{theorem}

\begin{proof}
The first inclusion is Theorem~\ref{teo:wf-accion-pseudodiferencial}. Now let $q\in\operatorname{Ell}(A)$ and suppose $q\notin\operatorname{WF}(A\mathbf{u})$. A microlocal parametrix gives $\mathbf{u}=B(A\mathbf{u})+R_{\mathbf{E}}\mathbf{u}$ near $q$. The first summand has no wavefront at $q$ by the preceding theorem, and neither does the second by \eqref{eq:microlocalmente-regularizante-wf}. Thus $q\notin\operatorname{WF}(\mathbf{u})$, which is the second inclusion. The remaining assertions reformulate these two inclusions.
\end{proof}

The wavefront set can also be defined entirely in the language of the pseudodifferential calculus.

\begin{theorem}[Characterization by operators of order zero]
\label{teo:caracterizacion-wf-psi0}
Let $M$ be a finite-dimensional smooth manifold without boundary, and let $\mathbf{E}\to M$ be a smooth complex vector bundle of finite rank. Let $\mathbf{u}\in\mathcal D'(M,\mathbf{E})$ and $q\in T^*M\setminus0_M$. The following are equivalent:
\begin{enumerate}[label=(\alph*)]
\item $q\notin\operatorname{WF}(\mathbf{u})$;
\item there exists $Q\in\Psi^0_{\mathrm{cl}}(M;\mathbf{E},\mathbf{E})$, properly supported and elliptic at $q$, such that $Q\mathbf{u}\in\Gamma(M,\mathbf{E})$.
\end{enumerate}
Consequently,
\begin{equation}
\label{eq:wf-interseccion-char-psi0}
\operatorname{WF}(\mathbf{u})
=
\bigcap
\bigl\{\operatorname{Char}(Q)
\mid Q\in\Psi^0_{\mathrm{cl}}(M;\mathbf{E},\mathbf{E}),\
Q\mathbf{u}\in\Gamma(M,\mathbf{E})\bigr\}.
\end{equation}
The intersection may be restricted to operators whose principal symbol is a scalar multiple of the identity of $\mathbf{E}$.
\end{theorem}

\begin{proof}
Suppose (a) holds. In a chart and local frame, choose a function $\varphi$ equal to one near the base point and a cone $\Gamma$ in which $\widehat{\varphi \mathbf{u}}$ decays rapidly. Let $\chi$ be a spatial cutoff equal to one on a smaller neighborhood, and $\vartheta$ a homogeneous cutoff of degree zero, angularly supported in $\Gamma$ and equal to one near the ray of $q$ at large frequencies. Quantize \[
q_0(x,\xi)=\chi(x)\vartheta(\xi)I_{\mathbf{E}}
\] and arrange proper support. In the region where $q_0$ contributes, the function $\vartheta\widehat{\varphi \mathbf{u}}$ decays rapidly. Fourier inversion and estimate \eqref{eq:fourier-localizacion-Au} show that $Q\mathbf{u}$ is smooth; the parts separated from the diagonal are also smooth. Moreover, the principal symbol of $Q$ is the identity near $q$, so $Q$ is elliptic there. This proves (b).

If (b) holds, microlocal elliptic regularity applied to $Q$ gives
\[
\operatorname{WF}(\mathbf{u})\cap\operatorname{Ell}(Q)
=\operatorname{WF}(Q\mathbf{u})\cap\operatorname{Ell}(Q)=\varnothing.
\]
Since $q\in\operatorname{Ell}(Q)$, we obtain (a). Formula \eqref{eq:wf-interseccion-char-psi0} is the set-theoretic expression of this equivalence. The construction in the first paragraph already uses a scalar symbol times $I_{\mathbf{E}}$, proving the last assertion.
\end{proof}

\subsection{The wavefront set of the kernel}

The conormal description of the kernel and the wavefront-set definition must agree on the diagonal. The following result identifies the possible directions exactly and explains, in microlocal terms, why a pseudodifferential operator relates a covector only to its negative in the input variable.

\begin{theorem}[Wavefront set of a pseudodifferential kernel]
\label{teo:wf-nucleo-pseudodiferencial}
Let $M$ be a finite-dimensional smooth manifold without boundary, and let $\mathbf{E},\mathbf{F}\to M$ be smooth complex vector bundles of finite rank. If $A\in\Psi^m_{\mathrm{cl,loc}}(M;\mathbf{E},\mathbf{F})$ and $\mathbf{K}_A$ is its Schwartz kernel, then
\begin{equation}
\label{eq:wf-kernel-conormal-diagonal}
\operatorname{WF}(\mathbf{K}_A)
\subseteq
N^*\Delta_M\setminus0_{M\times M}
=
\bigl\{(x,\xi;x,-\xi)\mid
(x,\xi)\in T^*M\setminus0_M\bigr\}.
\end{equation}
\end{theorem}

\begin{proof}
By Theorem~\ref{teo:equiv-pseudo-nucleo-conormal}, in a chart $\kappa\colon U\longrightarrow\Omega\subseteq\mathbb R^n$ and local frames of ranks $r_{\mathbf{E}},r_{\mathbf{F}}$, the kernel has the form, modulo a smooth function,
\[
K(x,y)
=(2\pi)^{-n}\operatorname{Os}\!\int_{\mathbb R_\xi^n}
e^{i(x-y)\cdot\xi}b(x,y,\xi)\,d\xi,
\qquad
b\colon\Omega\times\Omega\times\mathbb R_\xi^n
\longrightarrow
\operatorname{Hom}(\mathbb C^{r_{\mathbf{E}}},\mathbb C^{r_{\mathbf{F}}}),
\]
where $b$ is a classical amplitude of order $m$. We already know that $K$ is smooth if $x\neq y$. Now localize near a point of the diagonal with $\rho\in C_c^\infty$ and take the Fourier transform in $(x,y)$. Up to a constant depending only on $n$,
\begin{equation}
\label{eq:fourier-kernel-pseudo}
\widehat{\rho K}(\zeta,\eta)
=
\int_{\mathbb R^n}
\mathcal F_{(x,y)\to(\theta,\omega)}
\bigl(\rho b(\cdot,\cdot,\xi)\bigr)
(\zeta-\xi,\eta+\xi)\,d\xi.
\end{equation}
Integration by parts in $x$ and $y$ gives, for every $N$,
\[
\left\|
\mathcal F_{x,y}(\rho b)(\theta,\omega,\xi)
\right\|
\leq
C_N\langle(\theta,\omega)\rangle^{-N}
\langle\xi\rangle^m.
\]
Consider a closed cone of covectors $(\zeta,\eta)$ disjoint from the conormal bundle, that is, one in which
\[
|\zeta+\eta|
\geq c\bigl(|\zeta|+|\eta|\bigr).
\]
With $q=\displaystyle\frac{\zeta-\eta}{2}$ and $\xi=q+\tau$, we have the identity
\[
|\zeta-\xi|^2+|\eta+\xi|^2
=2|\tau|^2+\frac{1}{2}|\zeta+\eta|^2.
\]
In the cone under consideration, $|q|\leq C|\zeta+\eta|$. Inserting the preceding estimate into \eqref{eq:fourier-kernel-pseudo}, choosing $N$ arbitrarily large, and integrating in $\tau$ gives, for every $L$,
\[
|\widehat{\rho K}(\zeta,\eta)|
\leq C_L\langle(\zeta,\eta)\rangle^{-L}.
\]
By Definition~\ref{def:frente-onda-euclidiano}, no covector with $\zeta+\eta\neq0$ belongs to the wavefront set over the diagonal. The remaining covectors are exactly $(x,\xi;x,-\xi)$, in accordance with \eqref{eq:conormal-diagonal-identificacion}. Invariance under coordinate changes and multiplication by the smooth factors arising from frames and a density globalize the inclusion.
\end{proof}

If one adopts the convention of changing sign in the input covariable,
\[
\operatorname{WF}'(\mathbf{K}_A)
:=
\bigl\{(x,\xi;y,\eta)
\mid(x,\xi;y,-\eta)\in\operatorname{WF}(\mathbf{K}_A)\bigr\},
\]
the inclusion \eqref{eq:wf-kernel-conormal-diagonal} becomes
\[
\operatorname{WF}'(\mathbf{K}_A)
\subseteq
\bigl\{(x,\xi;x,\xi)\mid(x,\xi)\in T^*M\setminus0_M\bigr\}.
\]
Thus a pseudodifferential kernel transmits a singularity only at the same point and in the same cotangent direction. This assertion contains qualitative pseudolocality and, together with the microlocal parametrix, explains the elliptic equality \eqref{eq:igualdad-eliptica-wf}.

The symbolic construction of the parametrix and the conormal inclusion for the kernel's wavefront set correspond, under the translation $D=\partial$, to Theorems~18.1.9, 18.1.16, and 18.1.24 of \cite{HormanderIII}. These proofs make explicit the recurrence, control of remainders, and passage from Fourier analysis to the global formulation.

\section{Action on Sobolev spaces, uniform calculus, and Fredholm theory}
\label{sec:pseudo-sobolev-uniforme-fredholm}

The symbolic calculus acquires analytic content when translated into estimates between Sobolev spaces. In this section we write \[
H^s(M,\mathbf{E}):=H^{s,2}(M,\mathbf{E}),
\qquad s\in\mathbb R,
\] and use the bundle metrics and measure $d\lambda_{\mathbf{g}}$ to define the inner products on $L^2(M,\mathbf{E})$ and $L^2(M,\mathbf{F})$. In the complex case, we retain the previously fixed convention: the Hermitian inner product is linear in its first variable. We first prove the Euclidean order-zero estimate, then transfer it to manifolds, carefully distinguish the uniform noncompact case from the closed case, and deduce the Fredholm property.

\subsection{The Calderón--Vaillancourt estimate}

The decisive point is that an order-zero symbol defines a bounded operator between $L^2(\mathbb R^n,\mathbb C^r)$ and $L^2(\mathbb R^n,\mathbb C^s)$ even though its kernel need not be absolutely integrable. We give a proof using wave packets, because it exhibits finitely many seminorms and uses neither interpolation nor spectral theory.

\begin{theorem}[Calderón--Vaillancourt]
\label{pseudo:teo-calderon-vaillancourt}
\index{Calderon--Vaillancourt theorem@Calderón--Vaillancourt theorem}
Let $n,r,s\in\mathbb N$, let \[
a\in S^0_{1,0}\bigl(\mathbb R^n;
\operatorname{Hom}(\mathbb C^r,\mathbb C^s)\bigr)
\], and let $N_{\mathrm{CV}}:=4n+4$. There exists a constant $C=C(n,r,s)>0$ such that
\begin{equation}
\label{eq:calderon-vaillancourt-finito}
\|\operatorname{Op}(a)u\|_{L^2(\mathbb R^n,\mathbb C^s)}
\leq
Cq_{N_{\mathrm{CV}},0}(a)\|u\|_{L^2(\mathbb R^n,\mathbb C^r)}
\end{equation}
for every $u\in\mathcal S(\mathbb R^n,\mathbb C^r)$. Consequently, $\operatorname{Op}(a)$ extends uniquely to a continuous operator from $L^2(\mathbb R^n,\mathbb C^r)$ to $L^2(\mathbb R^n,\mathbb C^s)$.
\end{theorem}

\begin{proof}
It suffices to treat the scalar case: in finite rank, apply the argument to each matrix entry. Fix a nonzero $g\in C_c^\infty(\mathbb R^n)$ and define
\[
g_{z,\zeta}(x):=e^{ix\cdot\zeta}g(x-z),
\qquad
(Wu)(z,\zeta):=\langle u,g_{z,\zeta}\rangle_{L^2(\mathbb R^n)},
\qquad (z,\zeta)\in\mathbb R^{2n}.
\]
Plancherel's identity, applied first in $\zeta$ and then in $z$, gives Moyal's identity
\begin{equation}
\label{eq:moyal-paquetes-pseudo}
\|Wu\|_{L^2(\mathbb R^{2n})}^2
=(2\pi)^n\|g\|_{L^2(\mathbb R^n)}^2
\|u\|_{L^2(\mathbb R^n)}^2.
\end{equation}
By polarization, $W^*W=c_gI$, where $c_g=(2\pi)^n\|g\|_{L^2(\mathbb R^n)}^2$.

Consider the matrix kernel of $W\operatorname{Op}(a)W^*$:
\[
\mathcal M_a(z,\zeta;w,\eta)
:=
\langle\operatorname{Op}(a)g_{w,\eta},g_{z,\zeta}\rangle.
\]
The quantization formula of Definition~\ref{def:cuantizacion-kohn-nirenberg} writes it as an oscillatory integral in $(x,y,\theta)$ with phase \[
\Phi=x\cdot(\theta-\zeta)+y\cdot(\eta-\theta)
\] and amplitude
\[
(2\pi)^{-n}a(x,\theta)g(y-w)\overline{g(x-z)}.
\]
Set $L=n+1$. Since \[
(1-\Delta_x)^Le^{i\Phi}
=\langle\theta-\zeta\rangle^{2L}e^{i\Phi},
\qquad
(1-\Delta_y)^Le^{i\Phi}
=\langle\theta-\eta\rangle^{2L}e^{i\Phi},
\], integrate by parts $L$ times in $x$ and in $y$. Derivatives in $x$ hit $a$ and $g(x-z)$; derivatives in $y$ hit only $g(y-w)$. On the supports of these two factors,
\[
\langle z-w\rangle\leq C_g\langle x-y\rangle.
\]
To introduce this weight without changing the oscillatory integral, expand the polynomial $\langle z-w\rangle^{2L}=(1+|z-w|^2)^L$ using $z-w=(x-y)-(x-z)+(y-w)$. The factors in $x-z$ and $y-w$ are incorporated into the smooth compactly supported functions and their derivatives. The remaining factors are monomials $(x-y)^\gamma$ with $|\gamma|\leq2L$; use \[
(x-y)^\gamma e^{i\Phi}
=i^{-|\gamma|}D_\theta^\gamma e^{i\Phi}
\] to integrate by parts in $\theta$. Derivatives of the weights $\langle\theta-\zeta\rangle^{-2L}$ and $\langle\theta-\eta\rangle^{-2L}$ do not worsen their decay, and those hitting the symbol are controlled by $q_{4L,0}(a)$. This yields
\[
|\mathcal M_a(z,\zeta;w,\eta)|
\leq
C_gq_{4L,0}(a)\langle z-w\rangle^{-2L}
\int_{\mathbb R^n}
\langle\theta-\zeta\rangle^{-2L}
\langle\theta-\eta\rangle^{-2L}\,d\theta.
\]
Since $2L>n$, the elementary convolution inequality, together with Peetre's weight comparison in Lemma~\ref{lem:desigualdad-peetre-peso-japones}, gives
\[
\int_{\mathbb R^n}
\langle\theta-\zeta\rangle^{-2L}
\langle\theta-\eta\rangle^{-2L}\,d\theta
\leq C_L\langle\zeta-\eta\rangle^{-2L}.
\]
This inequality is proved by splitting the domain into the regions $\|\theta-\zeta\|\geq\displaystyle\frac{\|\zeta-\eta\|}{2}$ and $\|\theta-\eta\|\geq\displaystyle\frac{\|\zeta-\eta\|}{2}$. Thus
\begin{equation}
\label{eq:matriz-paquetes-cv}
|\mathcal M_a(z,\zeta;w,\eta)|
\leq
Cq_{4n+4,0}(a)
\langle z-w\rangle^{-2n-2}
\langle\zeta-\eta\rangle^{-2n-2}.
\end{equation}
The right-hand side is integrable in $(w,\eta)$ uniformly in $(z,\zeta)$, and also in $(z,\zeta)$ uniformly in $(w,\eta)$. Proposition~\ref{prop:criterio-schur-operadores-integrales}, applied with $X=Y=\mathbb R^{2n}$, Lebesgue measure, and $r=s=1$, shows that the integral operator with kernel $\mathcal M_a$ is bounded on $L^2(\mathbb R^{2n})$. Finally, \[
\operatorname{Op}(a)
=c_g^{-2}W^*\bigl(W\operatorname{Op}(a)W^*\bigr)W,
\] and \eqref{eq:moyal-paquetes-pseudo} yield \eqref{eq:calderon-vaillancourt-finito}. Density of $\mathcal S(\mathbb R^n,\mathbb C^r)$ in $L^2(\mathbb R^n,\mathbb C^r)$ gives the unique extension.
\end{proof}

The integer $4n+4$ is a sufficient choice using only finitely many derivatives. For other choices and proofs, see \cite{Taylor1974Pseudo,HormanderIII,Shubin2001}.

\subsection{Action on Sobolev scales}

Recall that \[
J^t=\operatorname{Op}(\langle\xi\rangle^t)
\colon H^t(\mathbb R^n)\longrightarrow L^2(\mathbb R^n)
\] is an isometric isomorphism with the unitary normalization fixed in Definition~\ref{def:transformada-fourier-L1} and Theorem~\ref{teo:plancherel-fourier}.

\begin{theorem}[Euclidean Sobolev mapping]
\label{pseudo:teo-mapeo-sobolev-euclidiano}
\index{pseudodifferential operator!action on Sobolev spaces}
Let $n,r,q\in\mathbb N$, $m,\sigma\in\mathbb R$, and \[
a\in S^m_{1,0}\bigl(
\mathbb R^n;\operatorname{Hom}(\mathbb C^r,\mathbb C^q)\bigr).
\]. Then
\[
\operatorname{Op}(a)\colon
H^{\sigma+m}(\mathbb R^n,\mathbb C^r)
\longrightarrow
H^\sigma(\mathbb R^n,\mathbb C^q)
\]
is continuous. Its norm is controlled by finitely many seminorms of $a$ in $S^m_{1,0}(\mathbb R^n;
\operatorname{Hom}(\mathbb C^r,\mathbb C^q))$.
\end{theorem}

\begin{proof}
The composition calculus of Theorem~\ref{teo:composicion-simbolica-pseudo} shows that \[
B:=J^\sigma\operatorname{Op}(a)J^{-(\sigma+m)}
\] is a pseudodifferential operator of order zero. The seminorms of a fixed number of terms in its symbol and of the remainder are controlled by finitely many seminorms of $a$. Theorem~\ref{pseudo:teo-calderon-vaillancourt} gives
\[
\|\operatorname{Op}(a)u\|_{H^\sigma(\mathbb R^n,\mathbb C^q)}
=\|J^\sigma\operatorname{Op}(a)u\|_{L^2(\mathbb R^n,\mathbb C^q)}
\leq C\|J^{\sigma+m}u\|_{L^2(\mathbb R^n,\mathbb C^r)}
=C\|u\|_{H^{\sigma+m}(\mathbb R^n,\mathbb C^r)}.
\]
\end{proof}

Before passing to manifolds, it is useful to isolate the other term that appears when localizing an operator: a smooth compactly supported kernel. Unlike the pseudodifferential part, this term improves regularity by arbitrarily many derivatives. The following estimate specifies how many kernel derivatives suffice to control each Sobolev mapping.

\begin{lemma}[Sobolev mapping of a smooth compactly supported kernel]
\label{pseudo:lem-mapeo-sobolev-nucleo-suave-euclidiano}
Let $U,V\subseteq\mathbb R^n$ be bounded open sets, and let $K\Subset U\times V$, $r,s\in\mathbb N$, and \[
k\in C^\infty
\bigl(U\times V;\operatorname{Hom}(\mathbb C^r,\mathbb C^s)\bigr),
\qquad
\operatorname{supp}k\subseteq K.
\]. For $\sigma,\tau\in\mathbb R$, define
\begin{equation}
\label{eq:numero-derivadas-nucleo-suave-euclidiano}
a_\sigma:=\max\{0,\lceil\sigma\rceil\},
\qquad
b_\tau:=\max\{0,\lceil-\tau\rceil\},
\qquad
N(n,\sigma,\tau):=a_\sigma+b_\tau.
\end{equation}
For $u\in H^\tau(V,\mathbb C^r)$, choose an extension $\widetilde u\in H^\tau(\mathbb R^n,\mathbb C^r)$ of $u$. For each $x\in U$, define $T_ku(x)\in\mathbb C^s$ by
\begin{equation}
\left\langle T_ku(x),w\right\rangle_{\mathbb C^s}
:=
\left\langle
\widetilde u,k(x,\,\cdot\,)^*w
\right\rangle_{
H^\tau(\mathbb R^n,\mathbb C^r),
H^{-\tau}(\mathbb R^n,\mathbb C^r)},
\qquad w\in\mathbb C^s,
\label{eq:definicion-dualidad-nucleo-suave-euclidiano}
\end{equation}
where $k(x,\,\cdot\,)^*w$ is extended by zero outside $V$ and the pairing is linear in the first variable, in accordance with this chapter's Hermitian convention. This definition is independent of the chosen extension and determines a continuous linear operator
\[
T_k\colon H^\tau(V,\mathbb C^r)
\longrightarrow H^\sigma(U,\mathbb C^s).
\]
The map $x\mapsto T_ku(x)$ is smooth and, for every multi-index $\alpha$,
\begin{equation}
\left\langle D_x^\alpha(T_ku)(x),w\right\rangle_{\mathbb C^s}
=
\left\langle
\widetilde u,\bigl(D_x^\alpha k\bigr)(x,\,\cdot\,)^*w
\right\rangle_{
H^\tau(\mathbb R^n,\mathbb C^r),
H^{-\tau}(\mathbb R^n,\mathbb C^r)}.
\label{eq:derivada-operador-nucleo-dualidad}
\end{equation}
If $u$ has a representative in $L^1_{\mathrm{loc}}(V,\mathbb C^r)$, then
\[
(T_ku)(x)=\int_Vk(x,y)u(y)\,dy,
\qquad x\in U.
\]
Moreover, there is a constant $C=C(n,\sigma,\tau,U,V,K,r,s)>0$ independent of $k$ such that
\begin{equation}
\label{eq:cota-mapeo-nucleo-suave-euclidiano}
\|T_k\|_{\mathcal L(
H^\tau(V,\mathbb C^r),H^\sigma(U,\mathbb C^s))}
\leq
C
\max_{|\alpha|+|\beta|\leq N(n,\sigma,\tau)}
\sup_{(x,y)\in K}
\|D_x^\alpha D_y^\beta k(x,y)\|.
\end{equation}
The choice of $N$ in \eqref{eq:numero-derivadas-nucleo-suave-euclidiano} is sufficient; smaller values may be valid but are not needed for later applications.
\end{lemma}

\begin{proof}
Spaces on open sets are understood with the quotient restriction norm of Definition~\ref{def:regularidad-intermedia-norma-r-n-abierto-s}, applied componentwise. Let $K_U:=\operatorname{pr}_1(K)\Subset U$ and $K_V:=\operatorname{pr}_2(K)\Subset V$. For $x\in U$ and $w\in\mathbb C^s$,
\[
k(x,\,\cdot\,)^*w\in C_c^\infty(V,\mathbb C^r),
\qquad
\operatorname{supp}\bigl(k(x,\,\cdot\,)^*w\bigr)\subseteq K_V.
\]
Its extension by zero belongs to $C_c^\infty(\mathbb R^n,\mathbb C^r)$ and, by Proposition~\ref{prop:regularidad-intermedia-propiedades-escala-hilbertiana}, to $H^{-\tau}(\mathbb R^n,\mathbb C^r)$. The right-hand side of \eqref{eq:definicion-dualidad-nucleo-suave-euclidiano} is therefore well defined. Since the Hermitian inner product is linear in its first variable, this expression is antilinear in $w$. If $(e_j)_{j=1}^s$ is an orthonormal basis of $\mathbb C^s$, the vector \[
\sum_{j=1}^s
\left\langle\widetilde u,k(x,\,\cdot\,)^*e_j\right\rangle_{
H^\tau(\mathbb R^n,\mathbb C^r),
H^{-\tau}(\mathbb R^n,\mathbb C^r)}e_j
\] satisfies \eqref{eq:definicion-dualidad-nucleo-suave-euclidiano} for every $w\in\mathbb C^s$. Nondegeneracy of the Hermitian inner product shows that it is the unique vector with this property.

Let us check independence of the extension. If $\widetilde u_1$ and $\widetilde u_2$ restrict to $u$ on $V$, then $d:=\widetilde u_1-\widetilde u_2$ vanishes as a distribution on $V$. Since $k(x,\,\cdot\,)^*w$ has compact support contained in $V$, the Sobolev pairing here agrees with the distributional action, and
\[
\left\langle d,k(x,\,\cdot\,)^*w\right\rangle_{
H^\tau(\mathbb R^n,\mathbb C^r),
H^{-\tau}(\mathbb R^n,\mathbb C^r)}=0.
\]
Thus \eqref{eq:definicion-dualidad-nucleo-suave-euclidiano} is independent of $\widetilde u$.

For each $w\in\mathbb C^s$, the map \[
x\longmapsto k(x,\,\cdot\,)^*w
\] is smooth from $U$ into the Fréchet space $C_{K_V}^\infty(\mathbb R^n,\mathbb C^r)$ of smooth functions supported in $K_V$. The inclusion of this space into $H^{-\tau}(\mathbb R^n,\mathbb C^r)$ is continuous. Composing with the pairing against $\widetilde u$ gives a smooth function of $x$ and allows differentiation within the duality. Taking $w$ successively in an orthonormal basis of $\mathbb C^s$ shows that $x\mapsto T_ku(x)$ is smooth and gives \eqref{eq:derivada-operador-nucleo-dualidad} for every multi-index $\alpha$.

If $u\in L^1_{\mathrm{loc}}(V,\mathbb C^r)$, pairing with a test function recovers the Lebesgue integral. By definition of the adjoint and the Hermitian convention,
\[
\begin{aligned}
\left\langle T_ku(x),w\right\rangle_{\mathbb C^s}
&=
\int_V
\left\langle u(y),k(x,y)^*w\right\rangle_{\mathbb C^r}\,dy\\
&=
\int_V
\left\langle k(x,y)u(y),w\right\rangle_{\mathbb C^s}\,dy.
\end{aligned}
\]
Since this holds for every $w\in\mathbb C^s$, the duality definition agrees with the usual integral asserted in the statement.

It remains to prove the estimate. Fix $\varepsilon>0$ and choose the extension so that
\[
\|\widetilde u\|_{H^\tau(\mathbb R^n,\mathbb C^r)}
\leq
\|u\|_{H^\tau(V,\mathbb C^r)}+\varepsilon.
\]
If $|\alpha|\leq a_\sigma$ and $\|w\|_{\mathbb C^s}=1$, duality, $b_\tau\geq-\tau$, and \eqref{eq:derivada-operador-nucleo-dualidad} give
\[
\begin{aligned}
\left|
\left\langle D_x^\alpha(T_ku)(x),w\right\rangle_{\mathbb C^s}
\right|
&\leq
\|\widetilde u\|_{H^\tau(\mathbb R^n,\mathbb C^r)}
\|\bigl(D_x^\alpha k\bigr)(x,\,\cdot\,)^*w\|_{H^{-\tau}(\mathbb R^n,\mathbb C^r)}\\
&\leq
\|\widetilde u\|_{H^\tau(\mathbb R^n,\mathbb C^r)}
\|\bigl(D_x^\alpha k\bigr)(x,\,\cdot\,)^*w\|_{H^{b_\tau}(\mathbb R^n,\mathbb C^r)}.
\end{aligned}
\]
For $q\in\mathbb N_0$, the polynomial comparison \[
\langle\xi\rangle^{2q}
\leq
C_{n,q}\sum_{|\beta|\leq q}|\xi^\beta|^2
\] and Plancherel's theorem~\ref{teo:plancherel-fourier}, applied componentwise, imply
\[
\|f\|_{H^q(\mathbb R^n,\mathbb C^r)}^2
\leq
C_{n,q}\sum_{|\beta|\leq q}
\|D^\beta f\|_{L^2(\mathbb R^n,\mathbb C^r)}^2.
\]
Since all preceding test functions are supported in $K_V$, this inequality with $q=b_\tau$ gives
\begin{equation}
\|\bigl(D_x^\alpha k\bigr)(x,\,\cdot\,)^*w\|_{
H^{b_\tau}(\mathbb R^n,\mathbb C^r)}
\leq
C_0
\max_{|\beta|\leq b_\tau}
\sup_{(z,y)\in K}
\|D_z^\alpha D_y^\beta k(z,y)\|.
\label{eq:cota-test-nucleo-por-derivadas}
\end{equation}
Here $C_0$ depends on $n,b_\tau,r,s$ and the volume of $K_V$, but not on $x,k,u$ or $w$.

Taking the supremum over unit vectors $w$, integrating in $x\in K_U$, summing over $|\alpha|\leq a_\sigma$, and applying the same Plancherel comparison with $q=a_\sigma$ yields
\[
\|T_ku\|_{H^{a_\sigma}(\mathbb R^n,\mathbb C^s)}
\leq
C_1
\|\widetilde u\|_{H^\tau(\mathbb R^n,\mathbb C^r)}
\max_{\substack{|\alpha|\leq a_\sigma\\|\beta|\leq b_\tau}}
\sup_{(x,y)\in K}
\|D_x^\alpha D_y^\beta k(x,y)\|.
\]
In this formula, extend $T_ku$ by zero outside $U$; the extension is smooth because $\operatorname{supp}(T_ku)\subseteq K_U\Subset U$. Since $a_\sigma\geq\sigma$, the inclusion $H^{a_\sigma}(\mathbb R^n,\mathbb C^s)
\hookrightarrow H^\sigma(\mathbb R^n,\mathbb C^s)$ and restriction to $U$ are continuous. Finally, $|\alpha|+|\beta|\leq a_\sigma+b_\tau=N(n,\sigma,\tau)$; let $\varepsilon\to 0^{+}$ to obtain \eqref{eq:cota-mapeo-nucleo-suave-euclidiano}.
\end{proof}

The local version requires a distinction between compact support and local regularity. For a compact set $K\Subset M$, set
\[
H^t_K(M,\mathbf{E}):=\{\mathbf{u}\in H^t_{\mathrm{loc}}(M,\mathbf{E})\mid\supp \mathbf{u}\subseteq K\},
\qquad
H^t_{\mathrm{comp}}(M,\mathbf{E}):=\bigcup_{K\Subset M}H^t_K(M,\mathbf{E}),
\]
with the usual inductive topology. The topology of $H^t_{\mathrm{loc}}(M,\mathbf{E})$ is defined by the seminorms $\|\chi \mathbf{u}\|_{H^t(M,\mathbf{E})}$, $\chi\in C_c^\infty(M)$, in accordance with Definition~\ref{def:regularidad-intermedia-sobolev-local-orden-real}.

\begin{proposition}[Local and compact-support versions]
\label{pseudo:prop-mapeos-sobolev-locales}
Let $(M,\mathbf{g})$ be a finite-dimensional smooth Riemannian manifold without boundary, not necessarily complete, and let $\mathbf{E},\mathbf{F}\to M$ be complex vector bundles of finite rank equipped with Hermitian bundle metrics and metric connections. Let $m,s\in\mathbb R$ and $A\in\Psi^m_{\mathrm{loc}}(M;\mathbf{E},\mathbf{F})$. Then
\[
A\colon H^{s+m}_{\mathrm{comp}}(M,\mathbf{E})
\longrightarrow H^s_{\mathrm{loc}}(M,\mathbf{F})
\]
is continuous. If $A$ is properly supported, there are also continuous extensions
\[
A\colon H^{s+m}_{\mathrm{comp}}(M,\mathbf{E})
\longrightarrow H^s_{\mathrm{comp}}(M,\mathbf{F}),
\qquad
A\colon H^{s+m}_{\mathrm{loc}}(M,\mathbf{E})
\longrightarrow H^s_{\mathrm{loc}}(M,\mathbf{F}).
\]
\end{proposition}

\begin{proof}
Let $\chi,\psi\in C_c^\infty(M)$. In finitely many charts and frames, $\chi A\psi$ is a matrix of Euclidean operators of order $m$ modulo a smooth kernel. Theorem~\ref{pseudo:teo-mapeo-sobolev-euclidiano} controls the pseudodifferential part, while Lemma~\ref{pseudo:lem-mapeo-sobolev-nucleo-suave-euclidiano} controls the smooth kernel and shows that it gains arbitrarily many derivatives. This proves the first mapping after choosing $\psi=1$ near the datum's support.

If $A$ is properly supported, the image of a compact set has support in another compact set. For the local mapping, given $\chi$, proper support allows $\psi\in C_c^\infty(M)$ to be chosen so that \[
\chi A=\chi A\psi
\] on a neighborhood of the support of $\chi$. The estimate already proved applies to $\psi u$. These estimates are precisely continuity for the stated inductive and projective topologies.
\end{proof}

\begin{corollary}[Local Sobolev elliptic regularity]
\label{pseudo:cor-regularidad-eliptica-sobolev-local}
Let $(M,\mathbf{g})$ be a finite-dimensional smooth Riemannian manifold without boundary, not necessarily complete, and let $\mathbf{E},\mathbf{F}\to M$ be complex vector bundles of finite rank equipped with Hermitian bundle metrics and metric connections. Let $m,s\in\mathbb R$, and let $A\in\Psi^m_{\mathrm{cl}}(M;\mathbf{E},\mathbf{F})$ be elliptic and properly supported. If $\mathbf{u}\in\mathcal D'(M,\mathbf{E})$ and $A\mathbf{u}\in H^s_{\mathrm{loc}}(M,\mathbf{F})$, then $\mathbf{u}\in H^{s+m}_{\mathrm{loc}}(M,\mathbf{E})$. In particular, for every $t\in\mathbb R$,
\[
\mathbf{u}\in H^t_{\mathrm{loc}}(M,\mathbf{E})
\quad\Longleftrightarrow\quad
A\mathbf{u}\in H^{t-m}_{\mathrm{loc}}(M,\mathbf{F}).
\]
\end{corollary}

\begin{proof}
Let $B\in\Psi^{-m}_{\mathrm{cl}}(M;\mathbf{F},\mathbf{E})$ be a parametrix. The identity $BA=I-R$, with $R\in\Psi^{-\infty}(M;\mathbf{E},\mathbf{E})$, gives $\mathbf{u}=B(A\mathbf{u})+R\mathbf{u}$. Proposition~\ref{pseudo:prop-mapeos-sobolev-locales} maps the first term into $H^{s+m}_{\mathrm{loc}}(M,\mathbf{E})$, and the second is smooth. This proves the first assertion. The remaining implication in the equivalence is the local mapping property of $A$. For differential operators, this recovers the regularity of Theorem~\ref{teo:regularidad-eliptica-local-operadores-haces}; the pseudodifferential calculus adds the localization in codirections of Theorem~\ref{teo:regularidad-eliptica-microlocal}.
\end{proof}

In particular, if $M$ is closed, every proper-support condition is automatic and a finite cover gives the following result.

\begin{corollary}[Mapping on a closed manifold]
\label{pseudo:cor-mapeo-sobolev-cerrada}
Let $(M,\mathbf{g})$ be a closed Riemannian manifold, and let $\mathbf{E},\mathbf{F}\to M$ be complex Hermitian bundles of finite rank equipped with metric connections. Let $m,s\in\mathbb R$ and $A\in\Psi^m(M;\mathbf{E},\mathbf{F})$. Then
\[
A_s\colon H^{s+m}(M,\mathbf{E})\longrightarrow H^s(M,\mathbf{F})
\]
is continuous. The norm of $A_s$ is controlled by finitely many symbol seminorms in a finite atlas and finitely many seminorms of the smoothing kernel.
\end{corollary}
\begin{proof}
Since $M$ is compact, the constant function $1$ belongs to $C_c^\infty(M)$. In the topology of $H^t_{\mathrm{loc}}(M,\mathbf{E})$ defined before Proposition~\ref{pseudo:prop-mapeos-sobolev-locales}, the seminorm corresponding to $\chi=1$ is
\[
 \|1\mathbf{u}\|_{H^t(M,\mathbf{E})}
 =\|\mathbf{u}\|_{H^t(M,\mathbf{E})}.
\]
Conversely, continuity of multiplication by smooth functions, obtained in coordinates using Lemma~\ref{lem:regularidad-intermedia-multiplicacion-afinada}, shows that each seminorm $\|\chi\mathbf{u}\|_{H^t(M,\mathbf{E})}$ is bounded by a multiple of $\|\mathbf{u}\|_{H^t(M,\mathbf{E})}$. Thus \[
 H^t_{\mathrm{loc}}(M,\mathbf{E})=H^t(M,\mathbf{E})
\] as topological spaces. Moreover, the support of every distribution on $M$ is a closed subset of a compact space. In the definition of $H^t_{\mathrm{comp}}(M,\mathbf{E})$ we may take $K=M$, giving \[
 H^t_{\mathrm{comp}}(M,\mathbf{E})=H^t(M,\mathbf{E})
\] also as topological spaces.

Proposition~\ref{pseudo:prop-mapeos-sobolev-locales}, with these two identifications and $t=s+m$, gives the continuous mapping
\[
 A_s\colon H^{s+m}(M,\mathbf{E})\longrightarrow H^s(M,\mathbf{F}).
\]

For the assertion on the norm, compactness of $M$ allows a finite atlas, unitary trivializations, and finitely many cutoffs to be chosen. In the proof of Proposition~\ref{pseudo:prop-mapeos-sobolev-locales}, each part localized near the diagonal is estimated by Theorem~\ref{pseudo:teo-mapeo-sobolev-euclidiano}; its constant depends on finitely many seminorms of the local symbol. The parts away from the diagonal have smooth kernel by Proposition~\ref{prop:equiv-simbolo-nucleo-regularizante}, and Lemma~\ref{pseudo:lem-mapeo-sobolev-nucleo-suave-euclidiano} controls them through finitely many seminorms of that kernel. Since only finitely many localizations occur, the norm of $A_s$ is controlled by the finite collection of seminorms in the statement.
\end{proof}

\subsection{Uniform pseudodifferential calculus with finite propagation}

The preceding results are local, and on a closed manifold a finite cover suffices to combine their constants. On a noncompact manifold, the same estimate may deteriorate as the chart changes. To prevent this, we fix a bounded-geometry localization system, require uniform symbol bounds, and also control propagation of the remainders.

On an infinite cover, requiring each local representative to be a symbol is insufficient: the constants in its estimates might grow with the chart. Without compactness there are different choices of global calculus. We work on manifolds without boundary and with bounded geometry of all orders. We consider the calculus with finite propagation and uniform local bounds, which suffices for later applications and is compatible with the model of \cite[Section~5]{EichhornSchmidFIO}. Chapter~\ref{cap:geometria-acotada} already constructed geodesic charts, quadratic systems, and synchronous frames with all derivatives controlled; Chapter~\ref{cap:trazas-geometria-acotada} proved independence of section norms from those choices. We develop only the new pseudodifferential information.

Fix an admissible geodesic quadratic system \[
\mathcal Q=(U_i,\phi_i,h_i,\mathbf{e}_i^{\mathbf{E}},\mathbf{e}_i^{\mathbf{F}})_{i\in I},
\qquad U_i=B_{\mathbf{g}}(p_i,2r),
\qquad \sum_{i\in I}h_i^2=1,
\] with the synchronous frames of Proposition~\ref{prop:marcos-sincronos-geodesicos-uniformes}. The radius $r$ and separation of the centers are those of Theorem~\ref{thm:trivializacion-geodesica}. For $\rho\geq0$, set
\[
I_{\mathcal Q}(i,\rho)
:=
\{j\in I\mid d_{\mathbf{g}}(U_i,U_j)\leq\rho\},
\qquad
N_{\mathcal Q}(\rho)
:=
\sup_{i\in I}\#I_{\mathcal Q}(i,\rho).
\]

\begin{lemma}[Uniform number of interactions at bounded distance]
\label{pseudo:lem-conteo-propagacion-uniforme}
Let $(M,\mathbf{g})$ be a Riemannian manifold without boundary of bounded geometry, and let $\mathcal Q$ be the preceding admissible geodesic quadratic system, of radius $r>0$ and finite overlap multiplicity. For every $\rho\geq0$, we have $N_{\mathcal Q}(\rho)<\infty$. More precisely, if \[
L_{\mathcal Q}
:=
\sup_{i\in I}\#\{j\in I\mid U_i\cap U_j\neq\varnothing\}
\] and \[
k(\rho):=\left\lceil\frac{\rho+4r}{r}\right\rceil+2,
\], then
\begin{equation}
\label{eq:conteo-propagacion-uniforme}
N_{\mathcal Q}(\rho)
\leq
1+L_{\mathcal Q}+\cdots+L_{\mathcal Q}^{k(\rho)}.
\end{equation}
\end{lemma}

\begin{proof}
If $d_{\mathbf{g}}(U_i,U_j)\leq\rho$, then $d_{\mathbf{g}}(p_i,p_j)\leq\rho+4r$. Completeness and Hopf--Rinow provide a minimizing geodesic between the centers. Divide it into at most $\left\lceil\displaystyle\frac{\rho+4r}{r}\right\rceil+1$ segments of length at most $r$. Each subdivision point belongs to some ball $B_{\mathbf{g}}(p_\nu,r)$ of the cover. Two centers chosen at consecutive subdivision points are at distance at most $3r$, so their balls of radius $2r$ intersect. Adding $i$ and $j$ at the endpoints gives a chain of at most $k(\rho)$ intersections joining $U_i$ to $U_j$. From any chart there are at most $L_{\mathcal Q}$ choices at each step. Counting all chains of lengths from zero to $k(\rho)$ gives \eqref{eq:conteo-propagacion-uniforme}. This count permits repetitions and therefore remains an upper bound.
\end{proof}

The preceding lemma ensures that, for a fixed propagation radius, each chart interacts with only a uniformly bounded number of charts. With this control, we can choose a diagonal representation of the operator and define seminorms that separately measure its symbolic part and smooth remainder.

Set $B:=B_{\mathrm{euc}}(0,2r)$. The functions in the geodesic system can be chosen with \(\operatorname{supp}(h_i\circ\phi_i^{-1})\subset
B_{\mathrm{euc}}\left(0,\displaystyle\frac{3r}{2}\right)\), as in the proof of Theorem~\ref{thm:trivializacion-geodesica}. If $\mathbf{u}\in\Gamma_c(\mathbf{E})$, let $\mathbf{u}_i\colon U_i\to\mathbb C^{r_{\mathbf{E}}}$ denote its column of components in the frame $\mathbf{e}_i^{\mathbf{E}}$, and define local analysis and synthesis by
\begin{align*}
\mathcal A_i^{\mathbf{E}} \mathbf{u}
&:=\widetilde{(h_i\mathbf{u}_i)\circ\phi_i^{-1}},\\
\mathcal S_i^{\mathbf{F}}v
&:=h_i\sum_{b=1}^{r_{\mathbf{F}}}(v^b\circ\phi_i)\mathbf{e}_{i,b}^{\mathbf{F}}.
\end{align*}
The tilde denotes extension by zero; this is smooth because the support of $h_i$ is uniformly separated from the boundary of $U_i$. The auxiliary cutoffs of \ref{item:B3} are inserted before extending an expression that does not already contain the factor $h_i$.

The locally finite sum \[
w_{\mathcal Q}(x,y):=\sum_{i\in I}h_i(x)h_i(y)
\] has uniformly bounded derivatives near the diagonal, and $w_{\mathcal Q}(x,x)=1$. Indeed, at most $L_{\mathcal Q}$ summands occur in each chart, and their derivatives are controlled by \ref{item:B2}. The mean value theorem in charts and uniform comparison of Euclidean and Riemannian distances give numbers $\delta_{\mathcal Q}>0$ and $c_{\mathcal Q}>0$ such that
\begin{equation}
\label{eq:peso-diagonal-sistema-uniforme}
w_{\mathcal Q}(x,y)\geq\frac{1}{2},
\qquad
d_{\mathbf{g}}(x,y)\leq2\delta_{\mathcal Q},
\end{equation}
and all derivatives of $w_{\mathcal Q}^{-1}$ are bounded in that region by constants independent of the chart. Decrease $\delta_{\mathcal Q}$ so that it is less than a quarter of the normal radius and the $2\delta_{\mathcal Q}$-neighborhood of each $\operatorname{supp}(h_i)$ remains in $U_i$.

Fix a function $\omega_{\mathcal Q}\in C^\infty(M\times M,[0,1])$ equal to one when $d_{\mathbf{g}}(x,y)\leq\displaystyle\frac{\delta_{\mathcal Q}}{3}$ and zero when $d_{\mathbf{g}}(x,y)\geq\displaystyle\frac{2\delta_{\mathcal Q}}{3}$. It is constructed by applying a fixed cutoff to $d_{\mathbf{g}}(x,y)^2$ within the normal radius; its covariant derivatives are uniformly bounded. In the chart $i$, write $\omega_i$ for its representative and set
\[
(\operatorname{Op}_{\omega_i}(a_i)v)(x)
:=
(2\pi)^{-n}\operatorname{Os}\!\iint_{B_y\times\mathbb R^n_\xi}
e^{i(x-y)\cdot\xi}\omega_i(x,y)a_i(x,\xi)v(y)\,dy\,d\xi,
\qquad x\in B.
\]
This formula uses the same chart and frames for $x$ and $y$. Define
\begin{equation}
\label{eq:pieza-diagonal-calculo-uniforme}
Q_i^{\mathcal Q}(a_i)
:=
\mathcal S_i^{\mathbf{F}}\operatorname{Op}_{\omega_i}(a_i)\mathcal A_i^{\mathbf{E}}.
\end{equation}
Its kernel is supported where $x,y\in\operatorname{supp}(h_i)$ and $d_{\mathbf{g}}(x,y)\leq\displaystyle\frac{2\delta_{\mathcal Q}}{3}$. In the usual abbreviated notation, \eqref{eq:pieza-diagonal-calculo-uniforme} is precisely \[
h_i\,\phi_i^*\operatorname{Op}_{\omega_i}(a_i)
(\phi_i^{-1})^*h_i,
\] with the identifications through $\mathbf{e}_i^{\mathbf{E}},\mathbf{e}_i^{\mathbf{F}}$ written in $\mathcal A_i^{\mathbf{E}}$ and $\mathcal S_i^{\mathbf{F}}$.

For a family $(a_i)_{i\in I}$, define the symbol seminorm of the representation by
\begin{equation}
\label{eq:seminorma-uniforme-simbolo}
q^{\mathrm u,\mathcal Q}_{N,m}((a_i))
:=
\sup_{i\in I}
\max_{|\alpha|+|\beta|\leq N}
\sup_{x\in B,\,\xi\in\mathbb R^n}
\langle\xi\rangle^{-m+|\beta|}
\|D_x^\alpha D_\xi^\beta a_i(x,\xi)\|.
\end{equation}
An operator with kernel $\mathbf{K}_R$ has propagation at most $\rho$ if
\[
 \operatorname{supp}\mathbf{K}_R
 \subseteq\{(x,y)\in M\times M\mid d_{\mathbf{g}}(x,y)\leq\rho\}.
\]
If $R$ has smooth kernel of propagation at most $\rho$, then
\[
\mathbf{K}_R\in C^\infty\!\left(
M_x\times M_y;\mathbf{F}_x\boxtimes \mathbf{E}_y^*
\right)
=C^\infty\!\left(
M_x\times M_y;\operatorname{Hom}(\mathbf{E}_y,\mathbf{F}_x)
\right).
\]
Equip this bundle with the product connection induced by the connection of $\mathbf{F}$ in the variable $x$ and the dual connection of $\mathbf{E}$ in the variable $y$. The expression $(\nabla_x)^p(\nabla_y)^q\mathbf{K}_R$ denotes iterated tensorial derivatives in each factor, with the new index first in its block; the two product connections commute. In the following expression, the norm is induced by $\mathbf{g}$, $\mathbf{h}_{\mathbf{F}}$, and $\mathbf{h}_{\mathbf{E}}$ on $(T_x^*M)^{\otimes p}\otimes(T_y^*M)^{\otimes q}\otimes
\operatorname{Hom}(\mathbf{E}_y,\mathbf{F}_x)$. Set
\begin{equation}
\label{eq:seminorma-uniforme-nucleo}
s_N^{\mathrm u}(R)
:=
\max_{p+q\leq N}
\sup_{(x,y)\in M\times M}
\| (\nabla_x)^p(\nabla_y)^q\mathbf{K}_R(x,y)\|_{\mathbf{g},\mathbf{h}_{\mathbf{F}},\mathbf{h}_{\mathbf{E}}}.
\end{equation}

\begin{definition}[Uniform diagonal representation]
\label{pseudo:def-calculo-uniforme}
\index{pseudodifferential operator!uniform}
\index{uniform pseudodifferential calculus@uniform pseudodifferential calculus}
Let $(M,\mathbf{g})$ be without boundary, and of bounded geometry, and let $\mathbf{E},\mathbf{F}\to M$ be complex vector bundles of finite rank equipped with Hermitian bundle metrics $\mathbf{h}_{\mathbf{E}},\mathbf{h}_{\mathbf{F}}$ and compatible connections of bounded geometry. A \textbf{uniform diagonal representation of order $m$ and residual radius $\rho$} of an operator $A\colon\Gamma_c(M,\mathbf{E})\to\Gamma(M,\mathbf{F})$ is a pair $\mathfrak r=((a_i)_{i\in I},R)$ such that
\begin{equation}
\label{eq:representacion-diagonal-calculo-uniforme}
A=A_{\mathrm{near}}+R,
\qquad
A_{\mathrm{near}}:=\sum_{i\in I}Q_i^{\mathcal Q}(a_i),
\end{equation}
$q^{\mathrm u,\mathcal Q}_{N,m}((a_i))<\infty$
for every $N$, and $R$ has propagation at most $\rho$ and satisfies $s_N^{\mathrm u}(R)<\infty$ for every $N$. Denote the set of these representations by $\operatorname{Rep}^{\mathcal Q}_{m,\rho}(A)$ and define their seminorms by
\begin{equation}
\label{eq:seminorma-representacion-uniforme}
P^{\mathcal Q}_{m,N}(\mathfrak r)
:=
q^{\mathrm u,\mathcal Q}_{N,m}((a_i))+s_N^{\mathrm u}(R).
\end{equation}
The operator belongs to $\Psi^m_{\mathrm u,\rho}(M;\mathbf{E},\mathbf{F})$ if this set is nonempty, and $\Psi^m_{\mathrm u}(M;\mathbf{E},\mathbf{F}):=\displaystyle\bigcup_{\rho\geq0}
\Psi^m_{\mathrm u,\rho}(M;\mathbf{E},\mathbf{F})$.
\end{definition}

If $\mathbf{E}=\mathbf{F}$, the identity belongs to the calculus with $a_i(x,\xi)=I$ and $R=0$, since $\operatorname{Op}_{\omega_i}(I)=I$ on the support of $h_i$ because $\omega_i=1$ on the diagonal, and
\[
\sum_{i\in I}Q_i^{\mathcal Q}(I)u
=
\sum_{i\in I}h_i^2u=u.
\]
This is the reason for using the same chart in both variables and a quadratic system.

The quantities \eqref{eq:seminorma-representacion-uniforme} depend on the representation. To obtain operator seminorms, fix $\rho\geq\delta_{\mathcal Q}$ and set
\begin{equation}
\label{eq:seminorma-intrinseca-operador-uniforme}
\|A\|^{\mathrm u,\mathcal Q}_{m,N;\rho}
:=
\inf_{\mathfrak r\in\operatorname{Rep}^{\mathcal Q}_{m,\rho}(A)}
P^{\mathcal Q}_{m,N}(\mathfrak r),
\end{equation}
with value $+\infty$ if no such representation exists. On $\Psi^m_{\mathrm u,\rho}(M;\mathbf{E},\mathbf{F})$ these are seminorms: addition and scalar multiplication are performed componentwise in representations, and choosing representations approximating the two infima gives the triangle inequality and homogeneity. They are the quotient seminorms of the linear map
\[
((a_i),R)\longmapsto\sum_{i\in I}Q_i^{\mathcal Q}(a_i)+R.
\]
The quotient seminorms separate points. If they all vanish for an operator $A$, then for each $j\geq1$ one can choose a representation of $A$ whose first $j$ seminorms sum to less than $2^{-j}$. On any fixed pair of charts, its kernels then converge to zero distributionally: for the pseudodifferential part this follows from continuity of quantization on test functions, and for the remainder from uniform convergence of the kernel. All representations have the same distributional kernel $\mathbf{K}_A$; by uniqueness, $\mathbf{K}_A=0$ and $A=0$. Thus each stage $\Psi^m_{\mathrm u,\rho}(M;\mathbf{E},\mathbf{F})$, with $\rho\geq\delta_{\mathcal Q}$, is a Hausdorff metrizable locally convex space, since the family of seminorms is indexed by $N\in\mathbb N_0$.

These stages are considered only as locally convex spaces; completeness, and hence the Fréchet property, would require proving that the chosen representation space is complete and that the kernel of the quotient map is closed in that realization. If $\rho_1\leq\rho_2$, inclusion of the first stage into the second is continuous, since a representation of radius $\rho_1$ is also one of radius $\rho_2$; its image need not be closed. In particular, the filtered union is not identified with a strict inductive limit or an LF-space.

The global structure used here is bornological and is fixed explicitly: a family $\mathscr B\subset\Psi^m_{\mathrm u}(M;\mathbf{E},\mathbf{F})$ is called \textbf{bounded in the uniform calculus} if there exists $\rho\geq\delta_{\mathcal Q}$ such that
$\mathscr B\subset\Psi^m_{\mathrm u,\rho}(M;\mathbf{E},\mathbf{F})$
and
\[
\sup_{A\in\mathscr B}
\|A\|^{\mathrm u,\mathcal Q}_{m,N;\rho}<\infty
\qquad\text{for every }N\in\mathbb N_0.
\]
This definition does not identify such families with the bounded sets of an inductive-limit topology. Only when $\mathcal Q$ and $\rho$ are fixed will we omit these superscripts from $q^{\mathrm u}$, $P$, and $\|\,\cdot\,\|^{\mathrm u}_{m,N;\rho}$.

Define $\Psi^{-\infty}_{\mathrm u}(M;\mathbf{E},\mathbf{F})$ to consist of finite-propagation operators $R\colon\Gamma_c(M,\mathbf{E})\to\Gamma(M,\mathbf{F})$ whose kernel \[
\mathbf{K}_R\in C^\infty
\bigl(M\times M;\operatorname{Hom}(\mathbf{E}_y,\mathbf{F}_x)\bigr)
\] satisfies \eqref{eq:seminorma-uniforme-nucleo} for every $N\in\mathbb N_0$. In the synchronous frames of the system $\mathcal Q$, this condition is equivalent to requiring all localized kernels $h_i(x)\mathbf{K}_R(x,y)h_j(y)$ to have partial derivatives of every order bounded by constants independent of $(i,j)$; the equivalence follows from uniform bounds on the connections, frame changes, and functions $h_i$. This is the required compatibility with the localization system. Finite propagation is global, not microlocal; by Hopf--Rinow, every operator in the calculus maps $\Gamma_c(M,\mathbf{E})$ into $\Gamma_c(M,\mathbf{F})$. If $M$ is noncompact, membership in $\Psi^{-\infty}_{\mathrm u}(M;\mathbf{E},\mathbf{F})$ does not imply compactness of the extension $L^2(M,\mathbf{E})\to L^2(M,\mathbf{F})$; Remark~\ref{pseudo:obs-regularizante-uniforme-no-compacto} gives a counterexample.

\begin{lemma}[Exact symbol of a uniformly smooth kernel]
\label{pseudo:lem-nucleos-uniformemente-suaves}
Let $r_{\mathbf{E}},r_{\mathbf{F}}\in\mathbb N$, let $B=B_{\mathrm{euc}}(0,2r)\subset\mathbb R^n$, and let \[
K_i\in C_c^\infty
\bigl(B\times B;
\operatorname{Hom}(\mathbb C^{r_{\mathbf{E}}},\mathbb C^{r_{\mathbf{F}}})\bigr),
\qquad i\in I,
\] be a family of kernels supported in a fixed compact subset of $B\times B$, with uniformly bounded derivatives. Their exact symbols
\[
r_i(x,\xi):=\int_{B_y} e^{-i(x-y)\cdot\xi}K_i(x,y)\,dy
\]
form a bounded family in $S^{-P}_{1,0}(B\times\mathbb R^n;
\operatorname{Hom}(\mathbb C^{r_{\mathbf{E}}},\mathbb C^{r_{\mathbf{F}}}))$ for each $P\in\mathbb R$. In particular, the definition of $\Psi^{-\infty}_{\mathrm u}(M;\mathbf{E},\mathbf{F})$ agrees with the intersection of all orders when a common radius is fixed and the preceding diagonal representations are used.
\end{lemma}

\begin{proof}
The identity avoiding spurious terms from differentiation of the phase is, for all multi-indices $\alpha,\beta$,
\begin{equation}
\label{eq:identidad-exacta-simbolo-nucleo-suave}
D_x^\alpha D_\xi^\beta r_i(x,\xi)
=
(-i)^{|\beta|}
\int_{B_y} e^{-i(x-y)\cdot\xi}
(D_x+D_y)^\alpha
\bigl[(x-y)^\beta K_i(x,y)\bigr]\,dy.
\end{equation}
To prove it, $D_\xi^\beta$ produces exactly $(-i)^{|\beta|}(x-y)^\beta$. Moreover, $(D_x+D_y)e^{-i(x-y)\cdot\xi}=0$; integrating the term $D_y$ by parts shows that differentiating the integral in $x$ is equivalent to applying $D_x+D_y$ to the remaining factor. Iterating this equality gives \eqref{eq:identidad-exacta-simbolo-nucleo-suave}.

Now use \[
\langle\xi\rangle^{2L}e^{-i(x-y)\cdot\xi}
=(1-\Delta_y)^Le^{-i(x-y)\cdot\xi}
\] and transfer $(1-\Delta_y)^L$ by integration by parts to the last factor of \eqref{eq:identidad-exacta-simbolo-nucleo-suave}. Since the supports remain in a fixed compact set, this yields
\[
\langle\xi\rangle^{2L}
\|D_x^\alpha D_\xi^\beta r_i(x,\xi)\|
\leq
C_{\alpha,\beta,L}
\max_{|\mu|+|\nu|\leq|\alpha|+|\beta|+2L}
\|D_x^\mu D_y^\nu K_i\|_{L^\infty(B\times B;\operatorname{Hom}(\mathbb C^{r_{\mathbf E}},\mathbb C^{r_{\mathbf F}}))},
\]
with a constant independent of $i$. Given $P$, choosing $2L\geq P+|\beta|$ proves the estimates for $S^{-P}_{1,0}(B\times\mathbb R^n;
\operatorname{Hom}(\mathbb C^{r_{\mathbf{E}}},\mathbb C^{r_{\mathbf{F}}}))$. Comparison of local partial and covariant derivatives, using the bounds on Christoffel symbols and connection matrices, identifies these estimates with \eqref{eq:seminorma-uniforme-nucleo}.

If $R\in\Psi^{-\infty}_{\mathrm u}(M;\mathbf{E},\mathbf{F})$, take $a_i=0$ in \eqref{eq:representacion-diagonal-calculo-uniforme}, so the operator belongs to every order with the same radius. For the converse, we do not assume that one representation works simultaneously at all orders. Let $A$ be a finite-propagation operator belonging to $\Psi^\mu_{\mathrm u}(M;\mathbf{E},\mathbf{F})$ for every $\mu\in\mathbb R$, and fix $p,q\in\mathbb N_0$. Choose $\mu<-n-p-q-2$ and a representation of $A$ of that order,
\[
A=\sum_{i\in I}Q_i^{\mathcal Q}(a_i^{(\mu)})+R_\mu.
\]
After applying $p$ derivatives in the output variable and $q$ in the input variable to the kernel of each diagonal piece, the oscillatory integral terms are bounded by $C_{p,q}\langle\xi\rangle^{\mu+p+q}$. This function is integrable in $\mathbb R^n$, and uniform multiplicity of the cover gives a chart-independent bound for the sum of pieces. The kernel of $R_\mu$ is already uniformly smooth. Consequently, the same distributional kernel $\mathbf{K}_A$ is represented by a kernel of class $C^{p,q}$ with the corresponding derivatives uniformly bounded. As $p,q$ varies, the resulting representations are compatible because they all represent the unique distributional kernel $\mathbf{K}_A$. Thus $\mathbf{K}_A$ is smooth and satisfies \eqref{eq:seminorma-uniforme-nucleo} for every $N$; together with finite propagation of $A$, this proves $A\in\Psi^{-\infty}_{\mathrm u}(M;\mathbf{E},\mathbf{F})$.
\end{proof}

The preceding representations control a family of symbols chart by chart. To construct parametrices, we must sum formal expansions using frequency cutoffs that can be chosen uniformly in the chart index. The following lemma is the uniform version of symbolic asymptotic summation.

\begin{lemma}[Uniform asymptotic summation]
\label{pseudo:lem-suma-asintotica-uniforme}
Let $r_{\mathbf{E}},r_{\mathbf{F}}\in\mathbb N$, let $\mathcal I$ be an index set, and let $m_0>m_1>\cdots$ be real numbers with $m_j\to-\infty$. In this lemma, abbreviate
\[
S^\mu
:=
S^\mu_{1,0}\bigl(
B\times\mathbb R^n;
\operatorname{Hom}(\mathbb C^{r_{\mathbf{E}}},\mathbb C^{r_{\mathbf{F}}})\bigr).
\]
For each $j\in\mathbb N_0$, let $(a^{(j)}_\nu)_{\nu\in\mathcal I}$ be a family in $S^{m_j}$ such that
\[
q^{\mathrm u}_{N,m_j}((a^{(j)}_\nu))<\infty
\qquad\text{for all }N,j.
\]
There is a family $(a_\nu)_{\nu\in\mathcal I}\subset S^{m_0}$, bounded in the seminorms of $S^{m_0}$, such that, for every $N\geq1$ and every $L\in\mathbb N_0$,
\begin{equation}
\label{eq:suma-asintotica-uniforme}
q^{\mathrm u}_{L,m_N}
\left((a_\nu)-\sum_{j=0}^{N-1}(a^{(j)}_\nu)\right)<\infty.
\end{equation}
If orders repeat, first group the families of equal order.
\end{lemma}

\begin{proof}
Take $\chi\in C^\infty(\mathbb R^n)$ vanishing for $\|\xi\|\leq1$ and equal to one for $\|\xi\|\geq2$. If $j>N$, the Leibniz rule and the support of the derivatives of $\chi(\varepsilon\xi)$ give
\begin{equation}
\label{eq:cota-corte-suma-uniforme}
q^{\mathrm u}_{L,m_N}
\bigl((\chi(\varepsilon\,\cdot)a^{(j)}_\nu)_\nu\bigr)
\leq
C_{L,N,j}\varepsilon^{m_N-m_j}
q^{\mathrm u}_{L,m_j}((a^{(j)}_\nu)),
\qquad 0<\varepsilon\leq1.
\end{equation}
Indeed, a term in which $D_\xi^\gamma$ hits the cutoff is supported in $\varepsilon^{-1}\leq\|\xi\|\leq2\varepsilon^{-1}$ if $|\gamma|>0$ and in $\|\xi\|\geq\varepsilon^{-1}$ if $\gamma=0$; in either case,
\[
\varepsilon^{|\gamma|}
\langle\xi\rangle^{m_j-m_N+|\gamma|}
\leq C_{N,j,\gamma}\varepsilon^{m_N-m_j}.
\]
The constant is uniform in $\nu$ because the seminorm in the last expression of \eqref{eq:cota-corte-suma-uniforme} is uniform.

Choose $\varepsilon_j\in(0,1]$ decreasing. Once $\varepsilon_0,\ldots,\varepsilon_{j-1}$ are fixed, the positive power of \eqref{eq:cota-corte-suma-uniforme} allows $\varepsilon_j$ to be chosen so that
\begin{equation}
\label{eq:eleccion-diagonal-cortes-uniformes}
q^{\mathrm u}_{L,m_N}
\bigl((\chi(\varepsilon_j\,\cdot)a^{(j)}_\nu)_\nu\bigr)
\leq2^{-j}
\end{equation}
for all pairs $0\leq N<j$ and $0\leq L\leq j$. Define
\[
a_\nu(x,\xi)
:=
\sum_{j=0}^\infty
\chi(\varepsilon_j\xi)a^{(j)}_\nu(x,\xi).
\]
For a fixed seminorm $q^{\mathrm u}_{L,m_0}$, the terms with $j>\displaystyle\max\{L,0\}$ form, by \eqref{eq:eleccion-diagonal-cortes-uniformes}, a series dominated by $\displaystyle \sum_{j=0}^{\infty}2^{-j}$; the remaining terms are finite in number. This proves convergence in all seminorms of $S^{m_0}$, uniformly in $\nu$.

Now fix $N\geq1$. We have the exact identity
\[
a_\nu-\sum_{j=0}^{N-1}a^{(j)}_\nu
=
\sum_{j=0}^{N-1}
(\chi(\varepsilon_j\xi)-1)a^{(j)}_\nu
+
\sum_{j=N}^\infty
\chi(\varepsilon_j\xi)a^{(j)}_\nu.
\]
The first sum has uniformly bounded support in $\xi$ and belongs to every symbol class. In the second, the terms with $N\leq j\leq\displaystyle\max\{N,L\}$ are finite in number and belong to $S^{m_N}$; the remainder converges in $q^{\mathrm u}_{L,m_N}$ by \eqref{eq:eleccion-diagonal-cortes-uniformes}. This proves \eqref{eq:suma-asintotica-uniforme} for each $L$.
\end{proof}

\begin{definition}[Uniform classicality and ellipticity]
\label{pseudo:def-elipticidad-uniforme-global}
\index{ellipticity!uniform}
Let $(M,\mathbf{g})$ be a Riemannian manifold without boundary of bounded geometry; let $\mathbf{E},\mathbf{F}\to M$ be complex vector bundles of finite rank equipped with Hermitian bundle metrics and compatible connections of bounded geometry; fix an admissible system $\mathcal Q$ as above, set $r_{\mathbf{E}}:=\operatorname{rk}\mathbf{E}$ and $r_{\mathbf{F}}:=\operatorname{rk}\mathbf{F}$, and let $m\in\mathbb R$. An operator $A\in\Psi^m_{\mathrm u}(M;\mathbf{E},\mathbf{F})$ is \textbf{uniformly classical} if it admits a diagonal representation whose symbols have expansions
\[
a_i\sim\sum_{j=0}^\infty a^{(i)}_{m-j}
\]
whose homogeneous components, restricted to $\|\xi\|=1$, have all derivatives bounded by constants independent of $i$, and if the remainders after $N$ components are uniformly bounded in
\[
S^{m-N}_{1,0}\bigl(
B\times\mathbb R^n;
\operatorname{Hom}(\mathbb C^{r_{\mathbf{E}}},\mathbb C^{r_{\mathbf{F}}})\bigr).
\]
We write $A\in\Psi^m_{\mathrm{u,cl}}(M;\mathbf{E},\mathbf{F})$.

Its principal symbol is the global section whose matrices in a chart are obtained from
\begin{equation}
\label{eq:simbolo-principal-representacion-diagonal}
\boldsymbol{\sigma}_m^\Psi(A)(x,\xi)
=
\sum_{i\in I}h_i(x)^2
\bigl[a^{(i)}_m(x,\xi)\bigr]_{\mathbf{e}_i^{\mathbf{E}},\mathbf{e}_i^{\mathbf{F}}},
\end{equation}
where each summand is transported to the corresponding fiber and covariable. The sum is finite at every point. Independence of the representation and system will be proved below by the exact transformation of kernels near the diagonal.

If $\mathbf{E}$ and $\mathbf{F}$ have the same rank, an operator $A\in\Psi^m_{\mathrm{u,cl}}(M;\mathbf{E},\mathbf{F})$ is \textbf{uniformly elliptic} if there exists $c>0$ such that
\begin{equation}
\bigl\|\boldsymbol{\sigma}_m^\Psi(A)(x,\xi)v\bigr\|_{\mathbf{h}_{\mathbf{F}}(x)}
\geq c\|v\|_{\mathbf{h}_{\mathbf{E}}(x)},
\qquad
x\in M,\quad \|\xi\|_{\mathbf{g}}=1,\quad v\in \mathbf{E}_x.
\label{eq:cota-elipticidad-uniforme-global}
\end{equation}
By homogeneity, for $\xi\neq0_x$ this is equivalent to
\[
\bigl\|\boldsymbol{\sigma}_m^\Psi(A)(x,\xi)^{-1}\bigr\|_{
\mathcal L((\mathbf{F}_x,\mathbf{h}_{\mathbf{F}}(x)),(\mathbf{E}_x,\mathbf{h}_{\mathbf{E}}(x)))}
\leq c^{-1}\|\xi\|_{\mathbf{g}}^{-m}.
\]
\end{definition}

The definition was formulated relative to an admissible geodesic system $\mathcal Q$. Before using the calculus, we must check that changing charts, frames, or cutoffs preserves the operator class and its bounded families. We therefore compare two admissible systems $\mathcal Q$ and $\mathcal Q^\sharp$, including their frames. The cross bounds in \ref{item:B1}, used in Lemma~\ref{lem:transferencia-uniforme-localizaciones} and Proposition~\ref{prop:independencia-espacios-haz}, provide the two finite cross multiplicities; let $L_\times$ denote their maximum. To specify how many derivatives are used, for $N\in\mathbb N_0$ and $\mu\in\mathbb R$ define
\[
L_N:=n+N+1,
\qquad
K_{N,\mu}:=n+N+1+\lceil|\mu|\rceil,
\qquad
\nu(N,\mu,n):=N+2L_N+2K_{N,\mu}+4.
\]
In the exact amplitude reduction formula \eqref{eq:simbolo-exacto-reduccion-amplitud}, $L_N$ integrations using $1-\Delta_\eta$ make the spatial variable integrable, and $K_{N,\mu}$ integrations using $1-\Delta_y$ make the frequency variable integrable after $N$ derivatives are taken. Thus the order-$N$ seminorm of the reduced symbol uses at most $\nu(N,\mu,n)$ derivatives of the amplitude. This number is not optimal, but it is fixed and independent of the chart.

Let $I^\sharp$ denote the index set of the system $\mathcal Q^\sharp$. Take a representation $\mathfrak r=((a_i),R)\in\operatorname{Rep}^{\mathcal Q}_{m,\rho}(A)$. Choose \(0<\delta<\displaystyle\frac{1}{4}\min\{\delta_{\mathcal Q},
\delta_{\mathcal Q^\sharp}\}\) and a uniform cutoff $\zeta(x,y)$ equal to one for $d_{\mathbf{g}}(x,y)\leq\delta$ and zero for $d_{\mathbf{g}}(x,y)\geq2\delta$. Split exactly as
\[
\mathbf{K}_{A_{\mathrm{near}}}
=
\zeta \mathbf{K}_{A_{\mathrm{near}}}
+(1-\zeta)\mathbf{K}_{A_{\mathrm{near}}}.
\]
The second summand is separated from the diagonal. The integrations by parts in the covariable used in the pseudolocality proof, with the uniform seminorms of $(a_i)$, show that it is a uniformly smooth kernel of propagation at most $\displaystyle\frac{2\delta_{\mathcal Q}}{3}$; incorporate it into $R$.

On the support of the first summand, $w_{\mathcal Q^\sharp}\geq\displaystyle\frac{1}{2}$ by \eqref{eq:peso-diagonal-sistema-uniforme}. Thus the identity
\begin{equation}
\label{eq:redistribucion-exacta-misma-carta}
\zeta \mathbf{K}_{A_{\mathrm{near}}}(x,y)
=
\sum_{a\in I^\sharp}
h_a^\sharp(x)
\frac{\zeta(x,y)\mathbf{K}_{A_{\mathrm{near}}}(x,y)}{w_{\mathcal Q^\sharp}(x,y)}
h_a^\sharp(y)
\end{equation}
redistributes the kernel exactly into pieces using the same chart $U_a^\sharp$ for $x$ and $y$. For fixed $a$, only the indices
\[
J(a):=\{i\in I\mid U_i\cap U_a^\sharp\neq\varnothing\},
\qquad \#J(a)\leq L_\times.
\]
contribute.

Write $\Phi_{ia}:=\phi_i\circ(\phi_a^\sharp)^{-1}$ on the overlap and
\[
B_{ia}(x,y):=\int_0^1D\Phi_{ia}(y+t(x-y))\,dt.
\]
Then $\Phi_{ia}(x)-\Phi_{ia}(y)=B_{ia}(x,y)(x-y)$, and, decreasing $\delta$ if necessary, $B_{ia}$ is invertible with uniformly bounded inverse. Let $G_{ia}^{\mathbf{E}}(y)$ and $G_{ia}^{\mathbf{F}}(x)$ be the frame-change matrices from the new frames to the old ones. After the change $\xi=B_{ia}(x,y)^T\eta$, the contribution of index $i$ to the term $a$ in \eqref{eq:redistribucion-exacta-misma-carta} has, before reduction, amplitude
\begin{equation}
\label{eq:amplitud-transferida-calculo-uniforme}
c_{ai}(x,y,\xi)
=\Theta_{ai}(x,y)
G_{ia}^{\mathbf{F}}(x)
a_i\bigl(\Phi_{ia}(x),B_{ia}(x,y)^{-T}\xi\bigr)
G_{ia}^{\mathbf{E}}(y)^{-1}
\frac{|\det D\Phi_{ia}(y)|}{|\det B_{ia}(x,y)|},
\end{equation}
where $\Theta_{ai}$ is the product of $\zeta$, the functions $h_i$ in both variables, the cutoff $\omega_i$, the necessary auxiliary cutoffs, and $w_{\mathcal Q^\sharp}^{-1}$. All these factors are expressed in the chart $a$ and have uniformly bounded derivatives. The formula separately includes the coordinate transition, both frame changes, and the Jacobian in the integration variable.

Exact reduction in Theorem~\ref{teo:amplitud-a-simbolo}, applied to $\displaystyle \sum_{i\in J(a)}c_{ai}$, produces a symbol $a_a^\sharp$ in the single chart $a$. The Leibniz rule, the Faà di Bruno formula of Theorem~\ref{faa di bruno multivariable}, the number of derivatives fixed by $\nu$, and $\#J(a)\leq L_\times$ give a representation $\mathfrak r^\sharp$ of $A$ for which
\begin{equation}
\label{eq:comparacion-seminormas-representacion-uniforme}
P^{\mathcal Q^\sharp}_{m,N}(\mathfrak r^\sharp)
\leq
C_{N,m}L_\times
P^{\mathcal Q}_{m,\nu(N,m,n)}(\mathfrak r).
\end{equation}
The constant depends on bounds up to order $\nu(N,m,n)+2$ for transitions and their inverses, frame changes and their inverses, cutoffs, bundle metrics, and densities, but not on $a$ or $i$. The new remainder has propagation at most $\rho^{\mathrm{com}}:=\displaystyle\max\{\rho,\delta_{\mathcal Q},
\delta_{\mathcal Q^\sharp}\}$.

Now take the infimum of \eqref{eq:comparacion-seminormas-representacion-uniforme} over all representations of radius at most $\rho$. This gives the intrinsic comparison
\begin{equation}
\label{eq:comparacion-seminormas-calculo-uniforme}
\|A\|^{\mathrm u,\mathcal Q^\sharp}_{m,N;\rho^{\mathrm{com}}}
\leq
C_{N,m}L_\times
\|A\|^{\mathrm u,\mathcal Q}_{m,\nu(N,m,n);\rho}.
\end{equation}
Interchanging $\mathcal Q$ and $\mathcal Q^\sharp$ gives the reverse inequality after using the common radius $\rho^{\mathrm{com}}$. This proves quantitative independence of the class, its filtration by radii, its bounded sets, and its locally convex topology. Applied to the principal homogeneous components of \eqref{eq:amplitud-transferida-calculo-uniforme}, the same formula gives the exact cotangent rule and shows that \eqref{eq:simbolo-principal-representacion-diagonal} is independent of the representation and system.

Having proved independence of choices, we must verify that the uniform calculus preserves the operations that motivated its construction. Write $A_h^*$ for the formal Hermitian adjoint of the pseudodifferential operator $A$, in accordance with the preceding local theorems. Reserve $T^\dagger$ for the Hilbert space adjoint of a bounded operator between two Hilbert spaces whose inner products are already fixed. For a representation of residual radius $\rho$, write
\[
\widehat\rho:=\max\left\{\rho,\frac{2\delta_{\mathcal Q}}{3}\right\};
\]
this bounds propagation of the full operator, including its diagonal part. Given $N\in\mathbb N_0$ and $\mu,\lambda\in\mathbb R$, define the derivative count and interaction factor
\begin{align}
\mathfrak n_N(\mu,\lambda)
&:=
\nu(N,\mu+\lambda,n)+\nu(N,\mu,n)+\nu(N,\lambda,n),
\label{eq:numero-derivadas-composicion-uniforme}\\
\mathfrak M_{\mathcal Q}(\rho,\sigma)
&:=
(1+L_{\mathcal Q})^2
\bigl(1+N_{\mathcal Q}(\rho+\sigma)\bigr)^2.
\label{eq:factor-interacciones-composicion-uniforme}
\end{align}
The second factor records both overlaps of diagonal pieces and bounded-distance interactions of residual kernels; applications require only that it be finite and depend on the stated data.

\begin{proposition}[Properties of the uniform calculus]
\label{pseudo:prop-calculo-uniforme}
Let $(M,\mathbf{g})$ be a Riemannian manifold without boundary of bounded geometry, and let $\mathbf{E},\mathbf{F},\mathbf{G}\to M$ be complex vector bundles of finite rank equipped with Hermitian bundle metrics and compatible connections of bounded geometry. Let $m,\ell\in\mathbb R$, $\rho_A,\rho_B\geq\delta_{\mathcal Q}$, $A\in\Psi^m_{\mathrm u,\rho_A}(M;\mathbf{E},\mathbf{F})$, and $B\in\Psi^\ell_{\mathrm u,\rho_B}(M;\mathbf{G},\mathbf{E})$. Then
\begin{enumerate}[label=(\roman*)]
\item $AB\in\Psi^{m+\ell}_{\mathrm u,\rho_{AB}}(M;\mathbf{G},\mathbf{F})$ with
\[
\rho_{AB}:=\rho_A+\rho_B;
\]
\item $A_h^*\in\Psi^m_{\mathrm u,\rho_A}(M;\mathbf{F},\mathbf{E})$;
\item the classes $\Psi^{-\infty}_{\mathrm u}(M;\mathbf{E},\mathbf{F})$ are stable under adjoints and form a two-sided ideal between the corresponding spaces within the finite-propagation calculus.
\end{enumerate}
In particular, composition is the continuous bilinear map of stages
\begin{equation}
\label{eq:composicion-etapas-calculo-uniforme}
\begin{split}
\Psi^{m_1}_{\mathrm u,\rho_1}(M;\mathbf{F},\mathbf{G})
\times
\Psi^{m_2}_{\mathrm u,\rho_2}(M;\mathbf{E},\mathbf{F})
&\longrightarrow
\Psi^{m_1+m_2}_{\mathrm u,\rho_1+\rho_2}(M;\mathbf{E},\mathbf{G}),\\
(A,B)&\longmapsto A\circ B,
\end{split}
\end{equation}
for $\rho_1,\rho_2\geq\delta_{\mathcal Q}$, and the adjoint is a continuous map from $\Psi^m_{\mathrm u,\rho}(M;\mathbf{E},\mathbf{F})$ to $\Psi^m_{\mathrm u,\rho}(M;\mathbf{F},\mathbf{E})$. For each $N$, there exists $C_{N,m,\ell}$ independent of local indices such that
\begin{equation}
\label{eq:cota-composicion-calculo-uniforme}
\|AB\|^{\mathrm u,\mathcal Q}_{m+\ell,N;\rho_{AB}}
\leq
C_{N,m,\ell}
\mathfrak M_{\mathcal Q}(\rho_A,\rho_B)
\|A\|^{\mathrm u,\mathcal Q}_{m,\mathfrak n_N(m,\ell);\rho_A}
\|B\|^{\mathrm u,\mathcal Q}_{\ell,\mathfrak n_N(m,\ell);\rho_B}.
\end{equation}
Moreover, with $\mathfrak a_N(m):=\nu(N,m,n)+N+2$,
\begin{equation}
\label{eq:cota-adjunto-calculo-uniforme}
\|A_h^*\|^{\mathrm u,\mathcal Q}_{m,N;\rho_A}
\leq
C_{N,m}\|A\|^{\mathrm u,\mathcal Q}_{m,\mathfrak a_N(m);\rho_A}.
\end{equation}
\end{proposition}

\begin{proof}
Take diagonal representations $\mathfrak r_A=((a_i),R_A)$ and $\mathfrak r_B=((b_j),R_B)$. Since $\rho_A,\rho_B\geq\delta_{\mathcal Q}$, we have $\widehat\rho_A=\rho_A$ and $\widehat\rho_B=\rho_B$. The kernel formula for composition is
\[
\mathbf{K}_{AB}(x,y)
=
\int_{M_z}\mathbf{K}_A(x,z)\mathbf{K}_B(z,y)\,d\lambda_{\mathbf{g}}(z).
\]
If the integrand is nonzero, then $d_{\mathbf{g}}(x,y)\leq\widehat\rho_A+\widehat\rho_B$, proving the propagation bound. Split as
\begin{equation}
\label{eq:cuatro-terminos-composicion-diagonal}
AB=A_{\mathrm{near}}B_{\mathrm{near}}
  +A_{\mathrm{near}}R_B+R_AB_{\mathrm{near}}+R_AR_B.
\end{equation}

Consider the first summand first. Its kernel is the sum of
\[
K_{ij}(x,y)
:=
\int_{M_z} \mathbf{K}_{Q_i(a_i)}(x,z)\mathbf{K}_{Q_j(b_j)}(z,y)\,d\lambda_{\mathbf{g}}(z).
\]
This occurs only if $U_i\cap U_j\neq\varnothing$; hence, for fixed $i$, there are at most $L_{\mathcal Q}$ values of $j$. Choose a cutoff $\zeta$ vanishing outside $d_{\mathbf{g}}(x,y)<\displaystyle\frac{\delta_{\mathcal Q}}{2}$ and equal to one near the diagonal. The part $(1-\zeta)K_{ij}$ is separated from the diagonal. Integration by parts in both covariables, before and after integration in $z$, gives all bounds in \eqref{eq:seminorma-uniforme-nucleo}; at most $L_{\mathcal Q}$ summands occur per chart. Incorporate this part into the remainder.

For $\zeta K_{ij}$, use the exact identity \eqref{eq:redistribucion-exacta-misma-carta} again:
\begin{equation}
\label{eq:redistribucion-composicion-misma-carta}
\zeta(x,y)\sum_{(i,j)\in I\times I}K_{ij}(x,y)
=
\sum_{a\in I} h_a(x)
\frac{\zeta(x,y)\sum_{(i,j)\in I\times I}K_{ij}(x,y)}{w_{\mathcal Q}(x,y)}h_a(y).
\end{equation}
After shrinking the support of $\zeta$ once, the choice of $\delta_{\mathcal Q}$ ensures that, if the summand $a,i,j$ is nonzero, all three points $x,z,y$ belong to $U_a$. Moreover, $U_i\cap U_a\neq\varnothing$ and $U_j\cap U_a\neq\varnothing$; for fixed $a$ there are at most $L_{\mathcal Q}^2$ pairs $(i,j)$.

Write $\Phi_{ia}=\phi_i\circ\phi_a^{-1}$ and
\[
B_{ia}(x,z)=\int_0^1
D\Phi_{ia}(z+t(x-z))\,dt.
\]
In the chart $a$, the change $\eta=B_{ia}(x,z)^{-T}\xi$ transforms the first piece into the amplitude
\begin{equation}
\label{eq:transporte-pieza-composicion-misma-carta}
c_{ia}(x,z,\xi)
=\Theta_{ia}(x,z)G_{ia}^{\mathbf{F}}(x)
a_i\bigl(\Phi_{ia}(x),B_{ia}(x,z)^{-T}\xi\bigr)
G_{ia}^{\mathbf{E}}(z)^{-1}
\frac{|\det D\Phi_{ia}(z)|}{|\det B_{ia}(x,z)|}.
\end{equation}
Here $\Theta_{ia}$ contains the corresponding cutoffs; for the second piece, with $\Phi_{ja}$, one obtains the same formula in the variables $(z,y)$ and with frame changes $G_{ja}^{\mathbf{E}},G_{ja}^{\mathbf{G}}$. Both intermediate matrices take values in the same frame $\mathbf{e}_a^{\mathbf{E}}$. This transformation takes place within the single chart $a$; it is not a quantization with base variables in different charts. Uniform bounds on transitions, frames, and Jacobians control all derivatives of both amplitudes. Exact reduction gives symbols $\widetilde a_{ia}(x,\xi)$ and $\widetilde b_{ja}(x,\xi)$, both defined in the chart $a$. Theorem~\ref{teo:composicion-simbolica-pseudo}, now applied in that single chart, yields
\[
\widetilde a_{ia}\#\widetilde b_{ja}
\sim
\sum_{\alpha\in\mathbb N_0^n}
\frac{i^{-|\alpha|}}{\alpha!}
D_\xi^\alpha\widetilde a_{ia}\,
D_x^\alpha\widetilde b_{ja}.
\]
Differentiating this formula up to order $N$, including reduction of both amplitudes, uses at most $\mathfrak n_N(m,\ell)$ derivatives of $(a_i)$ and $(b_j)$. Summing the at most $L_{\mathcal Q}^2$ pairs gives a symbol $c_a$ in the same chart $a$ and a uniformly smooth part.

The remaining three terms in \eqref{eq:cuatro-terminos-composicion-diagonal} are uniformly smooth. For $R_AR_B$, this is checked by differentiation under the integral; for each output chart, the integration domain is covered by at most $N_{\mathcal Q}(\widehat\rho_A+\widehat\rho_B)$ charts, each of uniformly bounded volume. In $A_{\mathrm{near}}R_B$, apply the pseudodifferential piece, in its single chart, to derivatives of the smooth kernel of $R_B$; in $R_AB_{\mathrm{near}}$, proceed in the other variable. The same count and Calderón--Vaillancourt after conjugation by powers of $\langle D\rangle$ give, for the constructed representation,
\begin{equation}
\label{eq:cota-representacion-composicion-uniforme}
P^{\mathcal Q}_{m+\ell,N}(\mathfrak r_{AB})
\leq
C_{N,m,\ell}\mathfrak M_{\mathcal Q}
 (\widehat\rho_A,\widehat\rho_B)
P^{\mathcal Q}_{m,\mathfrak n_N(m,\ell)}(\mathfrak r_A)
P^{\mathcal Q}_{\ell,\mathfrak n_N(m,\ell)}(\mathfrak r_B).
\end{equation}
Its remainder has radius at most $\rho_{AB}$. Choosing representations separately approximating the two infima and passing to those infima proves \eqref{eq:cota-composicion-calculo-uniforme}.

The adjoint of $Q_i^{\mathcal Q}(a_i)$ still uses the chart $i$ in both variables. Theorem~\ref{teo:simbolo-adjunto-pseudo} and exact reduction give a symbol $a_i^{(*)}$ with
\[
a_i^{(*)}(x,\xi)
\sim
\sum_{\alpha\in\mathbb N_0^n}
\frac{i^{-|\alpha|}}{\alpha!}
D_\xi^\alpha D_x^\alpha(a_i(x,\xi)^*).
\]
In nonunitary frames, insert the bundle metric matrices and the local density of $d\lambda_{\mathbf{g}}$. Their derivatives and those of their inverses are uniformly bounded. The difference between the exact adjoint kernel and the truncated quantization of $a_i^{(*)}$ is uniformly smooth and preserves propagation. This first proves the analogue of \eqref{eq:cota-adjunto-calculo-uniforme} for $P^{\mathcal Q}_{m,N}$ and then, taking the infimum, the asserted inequality. Finally, $\mathbf{K}_{R_h^*}(x,y)=\mathbf{K}_R(y,x)^*$ preserves the radius and all seminorms \eqref{eq:seminorma-uniforme-nucleo}. The arguments for the preceding three smooth terms, with one factor already smoothing, prove the two-sided ideal property.
\end{proof}

Stability under composition now leads to elliptic inversion. We first invert the principal symbol in each chart, correct lower-order terms by a uniform recurrence, and sum these corrections without losing control of propagation.

\begin{theorem}[Uniform parametrix]
\label{pseudo:teo-parametrix-uniforme}
Let $(M,\mathbf{g})$ be a Riemannian manifold without boundary of bounded geometry, and let $\mathbf{E},\mathbf{F}\to M$ be complex vector bundles of the same finite rank equipped with Hermitian bundle metrics and compatible connections of bounded geometry. If $m\in\mathbb R$ and $A\in\Psi^m_{\mathrm{u,cl}}(M;\mathbf{E},\mathbf{F})$ is uniformly elliptic, there exists $B\in\Psi^{-m}_{\mathrm u}(M;\mathbf{F},\mathbf{E})$ such that
\[
BA-I\in\Psi^{-\infty}_{\mathrm u}(M;\mathbf{E},\mathbf{E}),
\qquad
AB-I\in\Psi^{-\infty}_{\mathrm u}(M;\mathbf{F},\mathbf{F}).
\]
\end{theorem}

\begin{proof}
Fix a diagonal representation of $A$ of residual radius $\rho_A$ and use the convention $\widehat\rho_A$ introduced before Proposition~\ref{pseudo:prop-calculo-uniforme}. Localize the full symbol of $A$ in one chart at a time. More precisely, put $V:=B_{\mathrm{euc}}(0,\displaystyle\frac{7r}{4})\Subset B$ and use two
uniform spatial cutoffs supported in $B$ and equal to one on a neighborhood
of $\overline V$. Apply redistribution
\eqref{eq:redistribucion-exacta-misma-carta} and reduce the amplitude in chart
$i$. Let $a^{[i]}$ denote the full local symbol so obtained, to distinguish
it from the symbols of the individual diagonal pieces in a representation
of $A$. We use $a^{[i]}$ and its inverse only for $x\in V$. Its homogeneous families satisfy
\[
a^{[i]}(x,\xi)\sim\sum_{k=0}^\infty a^{[i]}_{m-k}(x,\xi),
\]
and, on $V$, $a_m^{[i]}$ is the matrix, in the frame $i$, of the global symbol $\boldsymbol{\sigma}_m^\Psi(A)$. Thus inversion is applied to the local matrix of the global principal symbol and uses uniform ellipticity directly.

For $j,N\in\mathbb N_0$, set
\[
A_{j,N}:=
\sup_{i\in I}\max_{|\gamma|\leq N}
\sup_{\substack{x\in V,\ \xi\in\mathbb R^n\\\|\xi\|=1}}
\|D^\gamma a^{[i]}_{m-j}(x,\xi)\|,
\qquad D^\gamma:=D_x^{\gamma_x}D_\xi^{\gamma_\xi}.
\]
Let $I_N$ be the same quantity for $(a_m^{[i]})^{-1}$. Differentiating $(a_m^{[i]})^{-1}a_m^{[i]}=I$, isolating the term $\eta=0$, and multiplying on the right by $(a_m^{[i]})^{-1}$ gives, for $|\gamma|\geq1$, the exact multi-index recurrence
\begin{equation}
\label{eq:recurrencia-derivadas-inversa-principal-uniforme}
D^\gamma(a_m^{[i]})^{-1}
=-
\sum_{0<\eta\leq\gamma}
\binom{\gamma}{\eta}
D^{\gamma-\eta}(a_m^{[i]})^{-1}
(D^\eta a_m^{[i]})(a_m^{[i]})^{-1}.
\end{equation}
In particular, $I_0\leq c^{-1}$ after uniform comparison of the fiber norms, and
\begin{equation}
\label{eq:cota-inductiva-derivadas-inversa-principal}
I_N
\leq
I_0+C_NI_0A_{0,N}\sum_{p=0}^{N-1}I_p.
\end{equation}
This recurrence proves inductively that $I_N<\infty$ for every $N$ and specifies which derivatives of the principal symbol occur.

The formal components of the left and right inverses are defined by
\begin{align}
b^{\mathrm L,(i)}_{-m}
&=b^{\mathrm R,(i)}_{-m}:=(a_m^{[i]})^{-1},
\label{eq:parametrix-uniforme-termino-inicial}\\
b^{\mathrm L,(i)}_{-m-j}
&:=-
\left[
\sum_{\substack{k,q\in\mathbb N_0,\ \alpha\in\mathbb N_0^n\\
                  k+q+|\alpha|=j,\ q<j}}
\frac{i^{-|\alpha|}}{\alpha!}
D_\xi^\alpha b^{\mathrm L,(i)}_{-m-q}
D_x^\alpha a^{[i]}_{m-k}
\right](a_m^{[i]})^{-1},
\quad j\geq1,
\label{eq:recurrencia-parametrix-uniforme-izquierda}\\
b^{\mathrm R,(i)}_{-m-j}
&:=-(a_m^{[i]})^{-1}
\left[
\sum_{\substack{k,q\in\mathbb N_0,\ \alpha\in\mathbb N_0^n\\
                  k+q+|\alpha|=j,\ q<j}}
\frac{i^{-|\alpha|}}{\alpha!}
D_\xi^\alpha a^{[i]}_{m-k}
D_x^\alpha b^{\mathrm R,(i)}_{-m-q}
\right].
\label{eq:recurrencia-parametrix-uniforme-derecha}
\end{align}
Here we have isolated $b_{-m-j}a_m$ in $b\#a$ and $a_mb_{-m-j}$ in $a\#b$, respectively; the other summands have $q<j$.

Define $B^{\mathrm L}_{j,N}$ and $B^{\mathrm R}_{j,N}$ as the supremum, over $i$, $x\in V$, $\|\xi\|=1$, and $|\gamma|\leq N$, of the derivatives of the corresponding components. The Leibniz rule and \eqref{eq:recurrencia-derivadas-inversa-principal-uniforme} give the quantified bounds
\begin{align}
B^{\mathrm L}_{j,N}
&\leq
C_{j,N}I_{N+j}
\sum_{\substack{k,q\in\mathbb N_0,\ \alpha\in\mathbb N_0^n\\
                  k+q+|\alpha|=j,\ q<j}}
B^{\mathrm L}_{q,N+j}A_{k,N+j},
\label{eq:induccion-cotas-parametrix-uniforme}\\
B^{\mathrm R}_{j,N}
&\leq
C_{j,N}I_{N+j}
\sum_{\substack{k,q\in\mathbb N_0,\ \alpha\in\mathbb N_0^n\\
                  k+q+|\alpha|=j,\ q<j}}
A_{k,N+j}B^{\mathrm R}_{q,N+j}.
\label{eq:induccion-cotas-parametrix-uniforme-derecha}
\end{align}
Induction on $j$ proves that both quantities are finite for every $j,N$. Each component has degree $-m-j$; homogeneity turns bounds on the sphere into uniform symbol estimates.

For the first Borel sum, fix radial cutoffs $\vartheta,\chi$ vanishing on $\|\xi\|\leq1$ and equal to one on $\|\xi\|\geq2$. Lemma~\ref{pseudo:lem-suma-asintotica-uniforme} chooses sequences $\varepsilon_j^{\mathrm L},\varepsilon_j^{\mathrm R}\to0^{+}$ for which
\begin{align}
\widetilde b_{\mathrm L}^{(i)}(x,\xi)
&:=\sum_{j=0}^\infty
\chi(\varepsilon_j^{\mathrm L}\xi)\vartheta(\xi)
b^{\mathrm L,(i)}_{-m-j}(x,\xi),
\label{eq:suma-borel-formal-parametrix-izquierda}\\
\widetilde b_{\mathrm R}^{(i)}(x,\xi)
&:=\sum_{j=0}^\infty
\chi(\varepsilon_j^{\mathrm R}\xi)\vartheta(\xi)
b^{\mathrm R,(i)}_{-m-j}(x,\xi)
\label{eq:suma-borel-formal-parametrix-derecha}
\end{align}
are families in $S^{-m}_{1,0}(V\times\mathbb R^n;
\operatorname{Hom}(\mathbb C^{r_{\mathbf{F}}},\mathbb C^{r_{\mathbf{E}}}))$ and, after subtracting the first $N$ summands, belong uniformly to the analogous class of order $-m-N$. These sums fix the expansion to be realized by the global inverse.
Choose a fixed $\eta\in C_c^\infty(V)$ equal to one on a neighborhood of
$\overline{B_{\mathrm{euc}}(0,\displaystyle\frac{3r}{2})}$. Multiply the sums by $\eta$ and extend
them by zero to $B$ before diagonal quantization. All symbol bounds remain
uniform. The factors $h_i$ in $Q_i^{\mathcal Q}$ see only the region where
$\eta=1$, so this extension changes none of the needed local inverse
identities. No inverse of the cutoff symbol is taken where the cutoffs vanish.

Their principal term defines, in the adopted diagonal model,
\begin{equation}
\label{eq:operador-inverso-principal-uniforme}
b_0^{(i)}(x,\xi):=\eta(x)\vartheta(\xi)(a_m^{[i]}(x,\xi))^{-1},
\qquad
B_0:=\sum_{i\in I}Q_i^{\mathcal Q}(b_0^{(i)}).
\end{equation}
Here $b_0^{(i)}$ is extended by zero from $V$ to $B$. In particular, diagonal quantization of $\widetilde b_{\mathrm L}^{(i)}$ or $\widetilde b_{\mathrm R}^{(i)}$ differs from $B_0$ at order $-m-1$. The equality $\displaystyle \sum_{i\in I} h_i^2=1$ shows that the principal symbol of $B_0$ is $\boldsymbol{\sigma}_m^\Psi(A)^{-1}$. Take $\rho_{B_0}:=\delta_{\mathcal Q}$ and
\[
R_0:=I-B_0A\in\Psi^{-1}_{\mathrm u,\varrho_0}(M;\mathbf{E},\mathbf{E}),
\qquad
S_0:=I-AB_0\in\Psi^{-1}_{\mathrm u,\varrho_0}(M;\mathbf{F},\mathbf{F}),
\qquad
\varrho_0:=\widehat\rho_A+\widehat\rho_{B_0}.
\]

We now specify the bounds on powers, including their radii. Let
\begin{equation}
\label{eq:coeficiente-recurrencia-potencias-parametrix}
\mathfrak C_N(\rho,\sigma;\mu,\lambda)
:=
C_{N,\mu,\lambda}\mathfrak M_{\mathcal Q}(\rho,\sigma),
\end{equation}
where $C_{N,\mu,\lambda}$ is the constant in \eqref{eq:cota-composicion-calculo-uniforme}. For $j\geq1$, set \begin{align*}
r_{j,N}&:=
\|R_0^j\|^{\mathrm u,\mathcal Q}_{-j,N;j\varrho_0},
&s_{j,N}&:=
\|S_0^j\|^{\mathrm u,\mathcal Q}_{-j,N;j\varrho_0},\\
\ell_{j,N}&:=
\|R_0^jB_0\|^{\mathrm u,\mathcal Q}_{-m-j,N;
j\varrho_0+\rho_{B_0}},
&d_{j,N}&:=
\|B_0S_0^j\|^{\mathrm u,\mathcal Q}_{-m-j,N;
\rho_{B_0}+j\varrho_0},
\end{align*} and $b_{0,N}:=\|B_0\|^{\mathrm u,\mathcal Q}_{-m,N;\rho_{B_0}}$, $a_N:=\|A\|^{\mathrm u,\mathcal Q}_{m,N;\rho_A}$. For $j=0$, adopt $\ell_{0,N}=d_{0,N}=b_{0,N}$. Exact cancellation of the principal symbols of $B_0A$ and $AB_0$, followed by amplitude reduction with one additional derivative, gives the base
\begin{equation}
\label{eq:base-recurrencia-defectos-parametrix-uniforme}
r_{1,N}+s_{1,N}
\leq
C_{N,m}\mathfrak M_{\mathcal Q}
(\widehat\rho_A,\widehat\rho_{B_0})
\left(1+
b_{0,\mathfrak n_{N+1}(-m,m)}
a_{\mathfrak n_{N+1}(-m,m)}\right).
\end{equation}
The composition construction followed by passage to the infimum gives, for $j\geq1$,
\begin{align}
r_{j+1,N}
&\leq
\mathfrak C_N(j\varrho_0,\varrho_0;-j,-1)
r_{j,\mathfrak n_N(-j,-1)}
r_{1,\mathfrak n_N(-j,-1)},
\label{eq:recurrencia-seminormas-potencias-defecto-izquierdo}\\
s_{j+1,N}
&\leq
\mathfrak C_N(j\varrho_0,\varrho_0;-j,-1)
s_{j,\mathfrak n_N(-j,-1)}
s_{1,\mathfrak n_N(-j,-1)},
\label{eq:recurrencia-seminormas-potencias-defecto-derecho}\\
\ell_{j,N}
&\leq
\mathfrak C_N(j\varrho_0,\rho_{B_0};-j,-m)
r_{j,\mathfrak n_N(-j,-m)}
b_{0,\mathfrak n_N(-j,-m)},
\label{eq:recurrencia-seminormas-terminos-neumann-izquierdos}\\
d_{j,N}
&\leq
\mathfrak C_N(\rho_{B_0},j\varrho_0;-m,-j)
b_{0,\mathfrak n_N(-m,-j)}
s_{j,\mathfrak n_N(-m,-j)}.
\label{eq:recurrencia-seminormas-terminos-neumann-derechos}
\end{align}
In particular, all four families of quantities are finite. The formulas simultaneously show growth of the radius, the increase $N\mapsto\mathfrak n_N$ in the derivative count, and the explicit appearance of $N_{\mathcal Q}(j\varrho_0+\varrho_0)$ within $\mathfrak M_{\mathcal Q}$. For $R_0^N$ and $S_0^N$, the first two recurrences, iterated from $r_{1,N}$ and $s_{1,N}$, are the required quantified bounds; in particular,
\begin{equation}
\label{eq:orden-radio-potencias-defectos-parametrix}
R_0^N\in\Psi^{-N}_{\mathrm u,N\varrho_0}(M;\mathbf{E},\mathbf{E}),
\qquad
S_0^N\in\Psi^{-N}_{\mathrm u,N\varrho_0}(M;\mathbf{F},\mathbf{F}).
\end{equation}

For $N\geq1$, introduce the finite sums
\begin{equation}
\label{eq:recurrencia-defectos-parametrix-uniforme}
C_N^{\mathrm L}:=\sum_{j=0}^{N-1}R_0^jB_0,
\qquad
C_N^{\mathrm R}:=B_0\sum_{j=0}^{N-1}S_0^j.
\end{equation}
The identities $B_0A=I-R_0$ and $AB_0=I-S_0$ give exactly
\begin{equation}
\label{eq:defectos-finitos-parametrix-uniforme}
I-C_N^{\mathrm L}A=R_0^N,
\qquad
I-AC_N^{\mathrm R}=S_0^N.
\end{equation}

Now use the representations produced recursively in the proof of Proposition~\ref{pseudo:prop-calculo-uniforme}. Let $(c_{j,i}^{\mathrm L})_i$ and $(c_{j,i}^{\mathrm R})_i$ be their symbol families near the diagonal for $R_0^jB_0$ and $B_0S_0^j$. Bounds \eqref{eq:recurrencia-seminormas-terminos-neumann-izquierdos} and \eqref{eq:recurrencia-seminormas-terminos-neumann-derechos}, before taking the infimum, verify the hypotheses of uniform asymptotic summation. With two further cutoff sequences, we obtain explicitly
\begin{align}
b_{\mathrm L}^{(i)}(x,\xi)
&:=\sum_{j=0}^\infty
\chi(\delta_j^{\mathrm L}\xi)c_{j,i}^{\mathrm L}(x,\xi),
&B_{\mathrm L}&:=\sum_{i\in I}Q_i^{\mathcal Q}(b_{\mathrm L}^{(i)}),
\label{eq:suma-borel-neumann-parametrix-izquierda}\\
b_{\mathrm R}^{(i)}(x,\xi)
&:=\sum_{j=0}^\infty
\chi(\delta_j^{\mathrm R}\xi)c_{j,i}^{\mathrm R}(x,\xi),
&B_{\mathrm R}&:=\sum_{i\in I}Q_i^{\mathcal Q}(b_{\mathrm R}^{(i)}).
\label{eq:suma-borel-neumann-parametrix-derecha}
\end{align}
Here only symbols of the parts near the diagonal are summed. The residual kernels of $R_0^jB_0$ and $B_0S_0^j$, whose radii grow with $j$, are not summed; this distinction keeps propagation fixed for $B_{\mathrm L}$ and $B_{\mathrm R}$.

For each $N$, the defining property of the two Borel sums and the finitely many remainders of terms with $j<N$ give
\begin{equation}
\label{eq:comparacion-sumas-finitas-borel-parametrix}
B_{\mathrm L}-C_N^{\mathrm L}
\in\Psi^{-m-N}_{\mathrm u}(M;\mathbf{F},\mathbf{E}),
\qquad
B_{\mathrm R}-C_N^{\mathrm R}
\in\Psi^{-m-N}_{\mathrm u}(M;\mathbf{F},\mathbf{E}).
\end{equation}
The radius of a representation of these differences may depend on $N$; the four operators $B_{\mathrm L}$, $B_{\mathrm R}$, $B_{\mathrm L}A-I$, and $AB_{\mathrm R}-I$, however, have fixed propagations. Combining \eqref{eq:comparacion-sumas-finitas-borel-parametrix} with \eqref{eq:defectos-finitos-parametrix-uniforme} gives, for every $N$,
\[
B_{\mathrm L}A-I\in\Psi^{-N}_{\mathrm u}(M;\mathbf{E},\mathbf{E}),
\qquad
AB_{\mathrm R}-I\in\Psi^{-N}_{\mathrm u}(M;\mathbf{F},\mathbf{F}).
\]
The two remainders are fixed finite-propagation operators belonging to $\Psi^{-N}_{\mathrm u}(M;\mathbf{E},\mathbf{E})$ and $\Psi^{-N}_{\mathrm u}(M;\mathbf{F},\mathbf{F})$, respectively, for every $N$. Apply the converse implication in Lemma~\ref{pseudo:lem-nucleos-uniformemente-suaves}: given $p,q$, choose $N>n+p+q+2$ and a representation of order $-N$ of the corresponding remainder. Absolute convergence of derivatives of its diagonal kernel and the uniformly smooth kernel of that representation provide the covariant bounds of orders $p,q$. Uniqueness of the distributional kernel makes these representations compatible as $p,q$ varies; hence the two remainders belong to $\Psi^{-\infty}_{\mathrm u}(M;\mathbf{E},\mathbf{E})$ and $\Psi^{-\infty}_{\mathrm u}(M;\mathbf{F},\mathbf{F})$, respectively.

If $B_{\mathrm L}A=I-R_{\mathbf{E}}$ and $AB_{\mathrm R}=I-R_{\mathbf{F}}$, then
\[
B_{\mathrm L}-B_{\mathrm R}
=B_{\mathrm L}R_{\mathbf{F}}-R_{\mathbf{E}}B_{\mathrm R}
\in\Psi^{-\infty}_{\mathrm u}(M;\mathbf{F},\mathbf{E}).
\]
Taking $B=B_{\mathrm L}$ and using the ideal property also on the right gives
\begin{equation}
\label{eq:parametrix-uniforme-residuos}
BA-I\in\Psi^{-\infty}_{\mathrm u}(M;\mathbf{E},\mathbf{E}),
\qquad
AB-I\in\Psi^{-\infty}_{\mathrm u}(M;\mathbf{F},\mathbf{F}).
\end{equation}
\end{proof}

The preceding parametrix remains a symbolic construction. To obtain global regularity, we must relate uniform seminorms to the Sobolev norms defined through the localization system. Lemma~\ref{pseudo:lem-mapeo-sobolev-nucleo-suave-euclidiano} controls the smooth remainders, while the diagonal part reduces to the Euclidean theorem.

\begin{theorem}[Uniform Sobolev mapping]
\label{pseudo:teo-mapeo-uniforme-sobolev}
Let $(M,\mathbf{g})$ be a Riemannian manifold without boundary of dimension $n$ and bounded geometry; let $\mathbf{E},\mathbf{F}\to M$ be complex vector bundles of finite rank equipped with Hermitian bundle metrics and compatible connections of bounded geometry, and fix the preceding admissible system $\mathcal Q$. If $m,s\in\mathbb R$, $\rho\geq\delta_{\mathcal Q}$, and $A\in\Psi^m_{\mathrm u,\rho}(M;\mathbf{E},\mathbf{F})$, then
\[
A\colon H^{s+m}(M,\mathbf{E})\longrightarrow H^s(M,\mathbf{F})
\]
is continuous. There exist an integer $N=N(n,m,s)$ and a constant $C_s$, depending on the uniform data of the system but not on a chart index or on $A$, such that
\begin{equation}
\label{eq:cota-mapeo-sobolev-uniforme-explicita}
\|A\|_{\mathcal L(H^{s+m}(M,\mathbf{E}),H^s(M,\mathbf{F}))}
\leq
C_s\bigl(1+N_{\mathcal Q}(\rho)\bigr)
\|A\|^{\mathrm u,\mathcal Q}_{m,N;\rho}.
\end{equation}
In particular, every $R\in\Psi^{-\infty}_{\mathrm u}(M;\mathbf{E},\mathbf{F})$ acts continuously from $H^t(M,\mathbf{E})$ to $H^s(M,\mathbf{F})$ for all $s,t\in\mathbb R$. All these estimates are uniform on bounded families in the calculus.
\end{theorem}

\begin{proof}
Let $\mathfrak r=((a_i),R)$ be a diagonal representation of $A$ of residual radius $\rho$. Introduce the analysis, diagonal block, and synthesis operators
\begin{align*}
\mathbf{A}_{\mathbf{E}}\mathbf{u}&:=(\mathcal A_i^{\mathbf{E}}\mathbf{u})_{i\in I},\\
\mathbf{T}_{(a_i)}(v_i)&:=
   (\operatorname{Op}_{\omega_i}(a_i)v_i)_{i\in I},\\
\mathbf{S}_{\mathbf{F}}(v_i)&:=\sum_{i\in I}\mathcal S_i^{\mathbf{F}}v_i.
\end{align*}
The part near the diagonal factors exactly as
\begin{equation}
\label{eq:factorizacion-analisis-sintesis-pseudodiferencial-uniforme}
A_{\mathrm{near}}
=\mathbf{S}_{\mathbf{F}}\mathbf{T}_{(a_i)}\mathbf{A}_{\mathbf{E}}.
\end{equation}
Analysis and synthesis are the retraction and coretraction operators of the preceding chapter: their continuity was established in \eqref{eq:sintesis-continua-cap13} and their equivalence with global norms in Lemma~\ref{lem:equivalencia-localizacion-cuadratica}. Reduction of the amplitude $\omega_i(x,y)a_i(x,\xi)$ uses finitely many derivatives, independently of $i$. Theorem~\ref{pseudo:teo-mapeo-sobolev-euclidiano}, applied componentwise, therefore provides $N_0=N_0(n,m,s)$ such that, for every $\mathbf{u}\in\Gamma_c(M,\mathbf{E})$,
\begin{equation}
\label{eq:cota-parte-diagonal-mapeo-uniforme}
\|A_{\mathrm{near}}\mathbf{u}\|_{H^s(M,\mathbf{F})}
\leq
C_s q^{\mathrm u,\mathcal Q}_{N_0,m}((a_i))
\|\mathbf{u}\|_{H^{s+m}(M,\mathbf{E})}.
\end{equation}

For the remainder, localize $h_iRh_j$. This can be nonzero only if $j\in I_{\mathcal Q}(i,\rho)$. Apply Lemma~\ref{pseudo:lem-mapeo-sobolev-nucleo-suave-euclidiano}, whose duality definition acts directly on restriction spaces, to the element $h_j\mathbf{u}\in H^{s+m}(U_j,\mathbf{E})$ on the fixed coordinate rectangles, with
\[
\sigma=s,\qquad \tau=s+m,\qquad
N_1
:=
\max\{0,\lceil s\rceil\}
+
\max\{0,\lceil-s-m\rceil\}.
\]
The supports of the localized kernels remain in a fixed compact subset of the rectangles, since the factors $h_i,h_j$ are uniformly separated from the chart boundaries. Uniform comparison of partial and covariant derivatives up to order $N_1$ gives
\[
\|h_iRh_j\mathbf{u}\|_{H^s(U_i,\mathbf{F})}
\leq
C_s s_{N_1}^{\mathrm u}(R)
\|h_j\mathbf{u}\|_{H^{s+m}(U_j,\mathbf{E})}.
\]
Each sum has at most $N_{\mathcal Q}(\rho)$ indices. Cauchy--Schwarz followed by reverse counting of the pairs $(i,j)$ yields
\begin{equation}
\label{eq:cota-residuo-mapeo-sobolev-uniforme}
\|R\mathbf{u}\|_{H^s(M,\mathbf{F})}
\leq
C_sN_{\mathcal Q}(\rho)s_{N_1}^{\mathrm u}(R)
\|\mathbf{u}\|_{H^{s+m}(M,\mathbf{E})}.
\end{equation}
With $N:=\displaystyle\max\{N_0,N_1\}$, \eqref{eq:cota-parte-diagonal-mapeo-uniforme} and \eqref{eq:cota-residuo-mapeo-sobolev-uniforme} first give
\[
\|A\|_{\mathcal L(H^{s+m}(M,\mathbf{E}),H^s(M,\mathbf{F}))}
\leq
C_s(1+N_{\mathcal Q}(\rho))P^{\mathcal Q}_{m,N}(\mathfrak r).
\]
Taking the infimum over all representations of radius $\rho$ proves \eqref{eq:cota-mapeo-sobolev-uniforme-explicita}. If $R$ is smoothing, use the representation $((0),R)$ and regard $R$ as an operator of order $t-s$; this gives the mapping from $H^t(M,\mathbf{E})$ to $H^s(M,\mathbf{F})$. Density of $\Gamma_c(M,\mathbf{E})$, proved using the same analysis--synthesis pair in the preceding chapter, completes the extensions.
\end{proof}

\begin{corollary}[Uniform elliptic estimate and regularity]
\label{pseudo:cor-regularidad-eliptica-uniforme-global}
Let $(M,\mathbf{g})$ be a Riemannian manifold without boundary of bounded geometry, and let $\mathbf{E},\mathbf{F}\to M$ be complex vector bundles of the same finite rank equipped with Hermitian bundle metrics and compatible connections of bounded geometry. Let $m\in\mathbb R$ and let $A\in\Psi^m_{\mathrm{u,cl}}(M;\mathbf{E},\mathbf{F})$ be uniformly elliptic. For $s,t\in\mathbb R$, there exists $C_{s,t}>0$ such that
\begin{equation}
\label{eq:estimacion-eliptica-uniforme-global}
\|\mathbf{u}\|_{H^{s+m}(M,\mathbf{E})}
\leq
C_{s,t}
\bigl(\|A\mathbf{u}\|_{H^s(M,\mathbf{F})}+\|\mathbf{u}\|_{H^t(M,\mathbf{E})}\bigr),
\qquad \mathbf{u}\in\Gamma_c(M,\mathbf{E}).
\end{equation}
If $\mathbf{u}\in H^t(M,\mathbf{E})$ and $A\mathbf{u}\in H^s(M,\mathbf{F})$ in the distributional sense, then $\mathbf{u}\in H^{s+m}(M,\mathbf{E})$.
\end{corollary}

\begin{proof}
Let $B$ be the parametrix of \eqref{eq:parametrix-uniforme-residuos} and write $BA=I-R_{\mathbf{E}}$. The preceding theorem gives
\[
\|\mathbf{u}\|_{H^{s+m}(M,\mathbf{E})}
\leq
\|BA\mathbf{u}\|_{H^{s+m}(M,\mathbf{E})}+\|R_{\mathbf{E}}\mathbf{u}\|_{H^{s+m}(M,\mathbf{E})}
\leq
C_{s,t}\bigl(\|A\mathbf{u}\|_{H^s(M,\mathbf{F})}+\|\mathbf{u}\|_{H^t(M,\mathbf{E})}\bigr),
\]
since $B$ has order $-m$ and $R_{\mathbf{E}}$ is uniformly smoothing. The same distributional identity proves the last assertion.
\end{proof}

\begin{remark}[Uniform does not mean compact]
\label{pseudo:obs-regularizante-uniforme-no-compacto}
\index{smoothing operator!uniform but noncompact}
On a noncompact manifold without boundary, an element of $\Psi^{-\infty}_{\mathrm u}(M;\mathbf{E},\mathbf{F})$ need not be compact or Fredholm. On $\mathbb R^n$, let $0\neq\kappa\in C_c^\infty(\mathbb R^n)$ and
\[
Ru:=\kappa*u.
\]
Its symbol $(2\pi)^{\displaystyle\frac{n}{2}}\widehat\kappa(\xi)$ belongs to $S^{-\infty}(\mathbb R^n;
\operatorname{Hom}(\mathbb C,\mathbb C))$ with the unitary normalization of Definition~\ref{def:transformada-fourier-L1}, and its kernel has finite propagation and translation-invariant bounds. Thus $R\in\Psi^{-\infty}_{\mathrm u}(\mathbb R^n;\mathbb C,\mathbb C)$. Choose $u\in C_c^\infty(\mathbb R^n)$ with $Ru\neq0$ and set $u_j(x)=u(x-je_1)$. Then $Ru_j$ is the corresponding translate of $Ru$; a subsequence with disjoint supports cannot converge. Thus $R$ is not compact. Moreover, $\widehat\kappa(\xi)\to0$ as $\|\xi\|\to\infty$; normalized functions with Fourier supports in disjoint balls escaping to infinity form an orthonormal sequence $(v_j)$ with $Rv_j\to0$. Since $v_j$ converges weakly to zero, its distance to any finite-dimensional subspace tends to one. This contradicts the lower bound for a Fredholm operator on the orthogonal complement of its kernel; thus $R$ is not Fredholm either.

For the same reason, existence of a parametrix with \emph{uniformly} smoothing remainders does not prove the Fredholm property. For example, $D_1=\partial_1\colon H^1(\mathbb R)\to L^2(\mathbb R)$ is uniformly elliptic but its range is not closed. Bounded geometry, even on a complete manifold without boundary, does not suffice for this conclusion: the preceding example requires additional control at infinity. We do not develop sufficient Fredholm conditions here.
\end{remark}

\subsection{Closed manifolds: compactness and Schatten classes}

On a noncompact manifold, a uniformly smooth kernel can translate to infinity without producing a compact operator. On a closed manifold this possibility disappears: smooth kernels are compact and, in fact, belong to Schatten classes. This difference will allow passage from a smoothing parametrix to the Fredholm property.

First fix the functional-analytic terminology used in this subsection. Let $H_1,H_2$ be separable Hilbert spaces and let $T\colon H_1\to H_2$ be compact. Its \textbf{singular values} are the eigenvalues, repeated according to multiplicity, arranged in nonincreasing order and padded with zeros when the rank is finite,
\[
s_j(T):=\lambda_j\bigl((T^\dagger T)^{\frac{1}{2}}\bigr),
\qquad j\geq1.
\]
For $1\leq p<\infty$, the \textbf{Schatten class} $\mathfrak S_p(H_1,H_2)$ consists of compact operators such that
\[
\|T\|_{\mathfrak S_p}
:=
\left(\sum_{j=1}^{\infty}s_j(T)^p\right)^{\frac{1}{p}}<\infty.
\]
In particular, $\mathfrak S_2$ is the Hilbert--Schmidt class and $\mathfrak S_1$ the trace class. If $(e_j)$ is an orthonormal basis of $H_1$, then
\begin{equation}
\label{eq:norma-hilbert-schmidt-base}
\|T\|_{\mathfrak S_2}^2
=
\sum_{j=1}^{\infty}\|Te_j\|_{H_2}^2;
\end{equation}
this quantity is independent of the basis, and its square root is
the \textbf{Hilbert--Schmidt norm}. If $H_1=H_2=H$ and $T\in\mathfrak S_1(H)$, its \textbf{trace} is
\begin{equation}
\label{eq:def-traza-schatten}
\operatorname{Tr}(T)
:=
\sum_{j=1}^{\infty}\langle Te_j,e_j\rangle_H.
\end{equation}
The series converges absolutely, and its sum is independent of the orthonormal basis. These assertions follow first for finite-rank operators by singular value decomposition and then by completeness in the norms $\mathfrak S_2$ and $\mathfrak S_1$.

\begin{proposition}[Square-integrable kernels and products of Schatten classes]
\label{pseudo:prop-nucleo-l2-schatten}
Let $(M,\mathbf{g})$ be a smooth closed Riemannian manifold of dimension $n$, and let $\mathbf{E},\mathbf{F}\longrightarrow M$ be complex vector bundles of finite rank equipped with Hermitian bundle metrics. An integral operator \[
(T_Ku)(x)=\int_{M_y} K(x,y)u(y)\,d\lambda_{\mathbf{g}}(y)
\] is Hilbert--Schmidt if and only if, up to almost-everywhere equality, its kernel belongs to
\[
L^2\bigl(M\times M;\operatorname{Hom}(\mathbf{E}_y,\mathbf{F}_x)\bigr).
\]
Here $\|K(x,y)\|_{\mathrm{HS}}$ is the Euclidean norm of the matrix of $K(x,y)\colon \mathbf{E}_y\to \mathbf{F}_x$ in orthonormal bases, independent of those bases. In this case,
\begin{equation}
\label{eq:norma-hs-nucleo-l2}
\|T_K\|_{\mathfrak S_2}^2
=
\int_{M_x\times M_y}\|K(x,y)\|_{\mathrm{HS}}^2
\,d\lambda_{\mathbf{g}}(x)\,d\lambda_{\mathbf{g}}(y).
\end{equation}
Moreover, if $C\in\mathfrak S_2(H_1,H_2)$ and $B\in\mathfrak S_2(H_2,H_3)$, then
\begin{equation}
\label{eq:producto-s2-s1}
BC\in\mathfrak S_1(H_1,H_3),
\qquad
\|BC\|_{\mathfrak S_1}
\leq
\|B\|_{\mathfrak S_2}\|C\|_{\mathfrak S_2}.
\end{equation}
\end{proposition}

\begin{proof}
The isometric identification between the Hilbert tensor product $L^2(M,\mathbf{F})\widehat\otimes\overline{L^2(M,\mathbf{E})}$ and Hilbert--Schmidt operators sends an elementary tensor $f\otimes\overline e$ to the operator $u\mapsto\langle u,e\rangle f$. In local orthonormal frames, this identification is precisely $K\mapsto T_K$; Parseval and Tonelli give \eqref{eq:norma-hs-nucleo-l2}. Surjectivity of the identification also proves the converse and uniqueness of the kernel in $L^2(M\times M;\operatorname{Hom}(\mathbf{E}_y,\mathbf{F}_x))$.

For the last assertion, use a singular value decomposition
\[
Cu=\sum_{j=1}^{\infty}s_j(C)\langle u,e_j\rangle f_j,
\]
where $(e_j)$ is orthonormal in $H_1$ and $(f_j)$ is an orthonormal system in $H_2$. Then
\[
BCu=\sum_{j=1}^{\infty}s_j(C)\langle u,e_j\rangle Bf_j
\]
is a nuclear representation. Indeed, a rank-one operator $u\mapsto\langle u,e\rangle h$ has $\mathfrak S_1$ norm equal to $\|e\|\|h\|$; by the triangle inequality, a series of such operators for which these quantities have finite sum converges in $\mathfrak S_1$. Now, by Cauchy--Schwarz and \eqref{eq:norma-hilbert-schmidt-base},
\[
\sum_{j=1}^{\infty}s_j(C)\|Bf_j\|
\leq
\left(\sum_{j=1}^{\infty}s_j(C)^2\right)^{\frac{1}{2}}
\left(\sum_{j=1}^{\infty}\|Bf_j\|^2\right)^{\frac{1}{2}}
\leq
\|C\|_{\mathfrak S_2}\|B\|_{\mathfrak S_2}.
\]
In the last inequality, $(f_j)$ was extended to an orthonormal basis of $H_2$. The partial sums therefore converge in trace-class norm and prove \eqref{eq:producto-s2-s1}.
\end{proof}

\begin{proposition}[Smoothing operators on a closed manifold]
\label{pseudo:prop-regularizantes-compactos}
\index{smoothing operator!compactness}
Let $(M,\mathbf{g})$ be a smooth closed Riemannian manifold, and let $\mathbf{E},\mathbf{F}\longrightarrow M$ be complex vector bundles of finite rank equipped with Hermitian bundle metrics and Hermitian connections. If $R\in\Psi^{-\infty}(M;\mathbf{E},\mathbf{F})$, then, for all $s,t\in\mathbb R$,
\[
R\colon H^s(M,\mathbf{E})\longrightarrow H^t(M,\mathbf{F})
\]
is continuous and compact. In fact, as an operator from $L^2(M,\mathbf{E})$ to $L^2(M,\mathbf{F})$, $R$ is trace class.
\end{proposition}

\begin{proof}
After localizing the kernel on coordinate rectangles, the duality definition of Lemma~\ref{pseudo:lem-mapeo-sobolev-nucleo-suave-euclidiano} applies directly to every element of $H^s(M,\mathbf{E})$ and allows differentiation within the pairing with the distribution. Thus $R$ maps $H^s(M,\mathbf{E})$ continuously into $H^N(M,\mathbf{F})$ for every $N$. Choose $N>t$. The inclusion $H^N(M,\mathbf{F})\hookrightarrow H^t(M,\mathbf{F})$ is compact by the Rellich--Kondrashov theorem, obtained from \ref{teo:rellich-orden-real} using a finite atlas; factorization proves compactness.

For the last assertion, localize the kernel on finitely many coordinate rectangles and extend each piece smoothly to a torus. Its double Fourier coefficients decay faster than any power. Thus the piece is a sum of rank-one operators with absolutely summable coefficients and is trace class. A finite sum preserves this property.
\end{proof}

\begin{lemma}[Diagonal formula for a continuous kernel]
\label{pseudo:lem-traza-nucleo-continuo}
Let $(M,\mathbf{g})$ be a smooth closed Riemannian manifold, and let $\mathbf{E}\longrightarrow M$ be a complex vector bundle of finite rank equipped with a Hermitian bundle metric. If $T\in\mathfrak S_1(L^2(M,\mathbf{E}))$ has a continuous integral kernel $\mathbf{K}_T(x,y)$, then
\[
\operatorname{Tr}(T)
=
\int_{M_x}\operatorname{tr}_{\mathbf{E}_x}\mathbf{K}_T(x,x)\,d\lambda_{\mathbf{g}}(x).
\]
\end{lemma}

\begin{proof}
Fix a finite cover $(V_i)$ with $\overline V_i$ contained in an open set $U_i$ equipped with a unitary frame $(\mathbf{e}_{i,a})$. Take finite measurable subdivisions $\mathcal P_\nu$ whose positive-measure cells lie in some $V_i$ and whose maximum diameter $\delta_\nu$ tends to zero; discard null cells. For each cell $Q\subset V_i$ and each vector of the local frame, let \[
\mathbf{v}_{Q,a}:=\lambda_{\mathbf{g}}(Q)^{-\frac{1}{2}}\mathbf1_Q\mathbf{e}_{i,a},
\], and denote by $P_\nu$ the orthogonal projection onto the space spanned by these vectors. If $\mathbf{u}\in\Gamma(M,\mathbf{E})$ and $\displaystyle \mathbf{u}=\displaystyle\sum_{a=1}^{r_{\mathbf E}} u_{i,a}\mathbf{e}_{i,a}$ in $U_i$, then, on $Q\subset V_i$,
\[
P_\nu \mathbf{u}
=
\sum_{a=1}^{r_{\mathbf E}}
\left(\frac{1}{\lambda_{\mathbf{g}}(Q)}
\int_{Q_y}u_{i,a}(y)\,d\lambda_{\mathbf{g}}(y)\right)\mathbf{e}_{i,a}.
\]
Uniformity of the moduli of continuity of the coefficients over this finite cover gives $\|P_\nu \mathbf{u}-\mathbf{u}\|_{L^2(M,\mathbf{E})}
\leq C\omega_{\mathbf{u}}(\delta_\nu)$. Density of smooth sections and $\|P_\nu\|\leq1$ imply that $P_\nu\to I$ strongly.

Since $T$ is trace class,
\[
\|T-P_\nu T P_\nu\|_{\mathfrak S_1}\longrightarrow0;
\]
indeed, $T-P_\nu TP_\nu=(I-P_\nu)T+P_\nu T(I-P_\nu)$, and both terms converge to zero in trace norm, first for $T$ of finite rank and then by approximating $T$ and $T^\dagger$ in $\mathfrak S_1$. Consequently,
\[
\operatorname{Tr}(T)
=
\lim_{\nu\to\infty}
\sum_{Q\in\mathcal P_\nu}\sum_{a=1}^{r_{\mathbf E}}
\frac{1}{\lambda_{\mathbf{g}}(Q)}
\int_{Q_x\times Q_y}
\big\langle \mathbf{K}_T(x,y)\mathbf{e}_{i,a}(y),\mathbf{e}_{i,a}(x)\big\rangle
\,d\lambda_{\mathbf{g}}(y)\,d\lambda_{\mathbf{g}}(x).
\]
On $\overline V_i\times\overline V_i$, set
\[
k_{i,a}(x,y)
:=
\big\langle \mathbf{K}_T(x,y)\mathbf{e}_{i,a}(y),\mathbf{e}_{i,a}(x)\big\rangle.
\]
These functions are uniformly continuous. If $\omega$ is the maximum of their moduli of continuity, the difference between the last sum and the diagonal integral is bounded by
\[
\operatorname{rk}(\mathbf{E})\lambda_{\mathbf{g}}(M)\omega(\delta_\nu),
\]
which tends to zero. Finally, $\displaystyle \sum_{a=1}^{r_{\mathbf E}} k_{i,a}(x,x)=\operatorname{tr}_{\mathbf{E}_x}\mathbf{K}_T(x,x)$, so the limit is independent of the frame and proves the formula.
\end{proof}

\begin{theorem}[Hilbert--Schmidt and trace-class thresholds]
\label{pseudo:teo-umbrales-schatten}
\index{Hilbert--Schmidt operator!pseudodifferential}
\index{trace-class operator!pseudodifferential}
Let $(M,\mathbf{g})$ be a smooth closed Riemannian manifold of dimension $n$, and let $\mathbf{E},\mathbf{F}\longrightarrow M$ be complex vector bundles of finite rank equipped with Hermitian bundle metrics and Hermitian connections. Let $m\in\mathbb R$.
\begin{enumerate}[label=(\roman*)]
\item If $A\in\Psi^m(M;\mathbf{E},\mathbf{F})$ and $m<-\displaystyle\frac{n}{2}$, then
\[
A\in\mathfrak S_2\bigl(L^2(M,\mathbf{E}),L^2(M,\mathbf{F})\bigr),
\]
that is, $A$ is Hilbert--Schmidt.
\item If $m<-n$, then
\[
A\in\mathfrak S_1\bigl(L^2(M,\mathbf{E}),L^2(M,\mathbf{F})\bigr).
\]
If $\mathbf{E}=\mathbf{F}$, its trace is given by the diagonal formula
\begin{equation}
\label{eq:traza-diagonal-pseudo}
\operatorname{Tr}(A)
=
\int_{M_x}\operatorname{tr}_{\mathbf{E}_x}\mathbf{K}_A(x,x)\,d\lambda_{\mathbf{g}}(x).
\end{equation}
\end{enumerate}
If, with the convention of Definition~\ref{def:cuantizacion-riemanniana-pseudo}, \[
A=\operatorname{Op}_{\mathbf{g},\nabla,\chi}(a)+R,
\qquad R\in\Psi^{-\infty}(M;\mathbf{E},\mathbf{E}),
\] and $m<-n$, then
\begin{equation}
\label{eq:traza-simbolo-geometrico-pseudo}
\operatorname{Tr}(A)
=
(2\pi)^{-n}
\int_{M_x}\int_{T_x^*M}\operatorname{tr}_{\mathbf{E}_x}a(x,\xi)
\,d\xi_x\,d\lambda_{\mathbf{g}}(x)
+
\int_{M_x}\operatorname{tr}_{\mathbf{E}_x}\mathbf{K}_R(x,x)\,d\lambda_{\mathbf{g}}(x).
\end{equation}
\end{theorem}

\begin{proof}
First localize a piece whose symbol $a(x,\xi)$ has compact support in $x$. Its kernel is
\[
K(x,y)
=(2\pi)^{-n}\int_{\mathbb R^n_\xi}
e^{i(x-y)\cdot\xi}a(x,\xi)\,d\xi.
\]
Plancherel in the variable $y$ gives
\[
\int_{\mathbb R^n_y}\|K(x,y)\|_{\mathrm{HS}}^2\,dy
=(2\pi)^{-n}
\int_{\mathbb R^n_\xi}\|a(x,\xi)\|_{\mathrm{HS}}^2\,d\xi.
\]
The right-hand side is integrable in $x$ when $2m<-n$. The Hilbert--Schmidt norm of an integral operator is the $L^2(M\times M;\operatorname{Hom}(\mathbf{E}_y,\mathbf{F}_x))$ norm of its kernel; a finite atlas and Proposition~\ref{pseudo:prop-regularizantes-compactos} for the remainders prove (i).

Suppose $m<-n$. Construct a scalar elliptic operator $Q_{\mathbf{E}}\in\Psi^{\displaystyle\frac{m}{2}}(M;\mathbf{E},\mathbf{E})$ with local full symbol $\langle\xi\rangle^{\displaystyle\frac{m}{2}}I_{\mathbf{E}}$ modulo lower-order terms and hence principal symbol $\|\xi\|^{\displaystyle\frac{m}{2}}I_{\mathbf{E}}$, and a parametrix $C_{\mathbf{E}}\in\Psi^{-\displaystyle\frac{m}{2}}(M;\mathbf{E},\mathbf{E})$. If $C_{\mathbf{E}}Q_{\mathbf{E}}=I-R_{\mathbf{E}}$, then
\[
A=(AC_{\mathbf{E}})Q_{\mathbf{E}}+AR_{\mathbf{E}}.
\]
The first two factors have order $\displaystyle\frac{m}{2}<-\displaystyle\frac{n}{2}$ and are Hilbert--Schmidt by (i); their product is trace class. The last term is smoothing and trace class by the preceding proposition. This proves (ii).

When $m<-n$, in each localized piece the integral defining the kernel has the uniform majorant $C\langle\xi\rangle^m\in L^1(\mathbb R^n)$. Dominated convergence shows that the kernel is jointly continuous in $(x,y)$; a finite sum of pieces and the smooth remainder preserve this property. Lemma~\ref{pseudo:lem-traza-nucleo-continuo}, applied to the trace-class operator already constructed, gives \eqref{eq:traza-diagonal-pseudo}. In geometric quantization, $j_x(0)=1$ and parallel transport on the diagonal is the identity; thus
\[
\mathbf{K}_{\operatorname{Op}_{\mathbf{g},\nabla,\chi}(a)}(x,x)
=(2\pi)^{-n}\int_{T_x^*M}a(x,\xi)\,d\xi_x.
\]
Adding the term $R$ proves \eqref{eq:traza-simbolo-geometrico-pseudo}.
\end{proof}

\subsection{Elliptic operators and Fredholm theory}

An elliptic parametrix gives inverse identities modulo smoothing operators. On a closed manifold these remainders are compact, so Atkinson's theory immediately turns symbolic inversion into a Fredholm statement. We begin with the form of this criterion suited to the two identities of a parametrix.

\begin{lemma}[Compact parametrix criterion]
\label{pseudo:lem-parametrix-compacta}
Let $X,Y$ be Hilbert spaces and let $T\in\mathcal L(X,Y)$ and $S\in\mathcal L(Y,X)$. If \[
ST-I_X\in\mathcal K(X),
\qquad
TS-I_Y\in\mathcal K(Y),
\], then $T$ is Fredholm in the sense of Definition~\ref{def:operador-fredholm-indice}.
\end{lemma}

\begin{proof}
Let $T^\dagger\colon Y\to X$ and $S^\dagger\colon X\to Y$ denote the bounded Hilbert space adjoints. These are not formal adjoints of differential operators: they depend on the specific inner products of $X$ and $Y$. If $Tu=0$, then $u=(I_X-ST)u$; the unit ball of $\ker T$ is therefore relatively compact. A normed space with relatively compact unit ball is finite-dimensional, so $\dim\ker T<\infty$. Applying the same argument to $T^\dagger$ and $S^\dagger$ gives $\dim\ker T^\dagger<\infty$.

Let us prove that the range is closed. On $(\ker T)^\perp$, if there were no $c>0$ with $\|Tu\|_Y\geq c\|u\|_X$, there would be a sequence $u_j\perp\ker T$, $\|u_j\|_X=1$, and $Tu_j\to0$. From \[
u_j=STu_j+(I_X-ST)u_j
\] and compactness of $I_X-ST$, we could extract a subsequence converging to $u$. Then $\|u\|=1$, $Tu=0$, and $u\perp\ker T$, a contradiction. The lower bound on $(\ker T)^\perp$ implies that $\operatorname{Ran}T$ is closed. Finally, \[
\operatorname{coker}T\cong(\operatorname{Ran}T)^\perp=\ker T^\dagger
\] is finite-dimensional.
\end{proof}

\begin{theorem}[Elliptic Fredholm theorem on a closed manifold]
\label{pseudo:teo-fredholm-eliptico}
\index{elliptic operator@elliptic operator!Fredholm property}
\index{Fredholm index@Fredholm index!elliptic operator}
Let $(M,\mathbf{g})$ be a smooth closed Riemannian manifold, and let $\mathbf{E},\mathbf{F}\longrightarrow M$ be complex vector bundles of finite rank equipped with Hermitian bundle metrics and Hermitian connections. Let $m,s\in\mathbb R$ and let $A\in\Psi^m_{\mathrm{cl}}(M;\mathbf{E},\mathbf{F})$ be elliptic. Then
\[
A_s\colon H^{s+m}(M,\mathbf{E})\longrightarrow H^s(M,\mathbf{F})
\]
is Fredholm. Its kernel and cokernel have smooth representatives:
\begin{align}
\ker A_s
&=
\{\mathbf{u}\in\Gamma(M,\mathbf{E})\mid A\mathbf{u}=0\},
\label{eq:nucleo-suave-pseudo}\\
\bigl(\operatorname{coker}A_s\bigr)'
&\cong
\{\mathbf{w}\in\Gamma(M,\overline{\mathbf{F}})\mid
  \overline{A_h^*}\,\mathbf{w}=0\}
\cong
\overline{\{\mathbf{v}\in\Gamma(M,\mathbf{F})\mid A_h^*\mathbf{v}=0\}}.
\label{eq:conucleo-suave-pseudo}
\end{align}
In particular,
\begin{equation}
\label{eq:indice-nucleos-suaves-pseudo}
\operatorname{ind}(A_s)
=
\dim\ker\bigl(A\colon\Gamma(M,\mathbf{E})\longrightarrow\Gamma(M,\mathbf{F})\bigr)
-
\dim\ker\bigl(A_h^*\colon\Gamma(M,\mathbf{F})\longrightarrow\Gamma(M,\mathbf{E})\bigr).
\end{equation}
\end{theorem}

\begin{proof}
The elliptic construction provides a parametrix $B\in\Psi^{-m}_{\mathrm{cl}}(M;\mathbf{F},\mathbf{E})$ and remainders
\[
BA=I-R_{\mathbf{E}},
\qquad
AB=I-R_{\mathbf{F}},
\qquad
R_{\mathbf{E}}\in\Psi^{-\infty}(M;\mathbf{E},\mathbf{E}),
\qquad
R_{\mathbf{F}}\in\Psi^{-\infty}(M;\mathbf{F},\mathbf{F}).
\]
By Corollary~\ref{pseudo:cor-mapeo-sobolev-cerrada}, $B\colon H^s(M,\mathbf{F})\longrightarrow H^{s+m}(M,\mathbf{E})$ is continuous, and by Proposition~\ref{pseudo:prop-regularizantes-compactos} both remainders are compact on the corresponding spaces. Lemma~\ref{pseudo:lem-parametrix-compacta} proves that $A_s$ is Fredholm.

If $\mathbf{u}$ is a distribution and $A\mathbf{u}=0$, the first identity gives $\mathbf{u}=R_E\mathbf{u}$. Since $R_{\mathbf{E}}$ is smoothing, $\mathbf{u}$ is smooth. This proves \eqref{eq:nucleo-suave-pseudo} and shows that the kernel is independent of the Sobolev order. By Proposition~\ref{prop:dualidad-sobolev-haces-geometria-acotada}, applied on the closed manifold, the complex-linear dual of $H^s(M,\mathbf{F})$ is $H^{-s}(M,\overline{\mathbf{F}})$. Define the conjugate of the Hermitian adjoint by
\[
 \overline{A_h^*}\colon\Gamma(M,\overline{\mathbf{F}})\longrightarrow
 \Gamma(M,\overline{\mathbf{E}}),
 \qquad
 \overline{A_h^*}(\overline{\mathbf{v}}):=\overline{A_h^*\mathbf{v}}.
\]
It is complex-linear. In charts it is obtained by conjugating the coefficients and kernel of $A_h^*$; hence it belongs to $\Psi^m_{\mathrm{cl}}(M;\overline{\mathbf{F}},\overline{\mathbf{E}})$ and is elliptic. The defining identity of the Hermitian adjoint becomes the bilinear identity
\[
 \langle A\mathbf{u},\mathbf{w}\rangle_{\mathbf{F},\overline{\mathbf{F}}}
 =\langle \mathbf{u},\overline{A_h^*}\mathbf{w}\rangle_{\mathbf{E},\overline{\mathbf{E}}}.
\]
Consequently, the annihilator of the range of $A_s$ is identified with
\[
\ker\left(
\overline{A_h^*}\colon H^{-s}(M,\overline{\mathbf{F}})
\longrightarrow H^{-s-m}(M,\overline{\mathbf{E}})
\right).
\]
The conjugate operator is elliptic, and the same smoothing argument makes its solutions smooth. Pointwise conjugation identifies its kernel with the conjugate space of $\ker A_h^*$. This gives \eqref{eq:conucleo-suave-pseudo} and, together with the definition of the index, yields \eqref{eq:indice-nucleos-suaves-pseudo}.
\end{proof}

\begin{corollary}[Independence of the Sobolev exponent]
\label{pseudo:cor-independencia-s-indice}
Let $(M,\mathbf{g})$ be a smooth closed Riemannian manifold; let $\mathbf{E},\mathbf{F}\longrightarrow M$ be complex vector bundles of finite rank equipped with Hermitian bundle metrics and Hermitian connections; and let $m\in\mathbb R$ and let $A\in\Psi^m_{\mathrm{cl}}(M;\mathbf{E},\mathbf{F})$ be elliptic. Then \[
\operatorname{ind}\left(
A\colon H^{s+m}(M,\mathbf{E})\longrightarrow H^s(M,\mathbf{F})
\right)
\] is independent of $s\in\mathbb R$.
\end{corollary}

\begin{proof}
The two spaces of smooth solutions appearing in \eqref{eq:indice-nucleos-suaves-pseudo} are independent of $s$.
\end{proof}

\subsection{Analytic stability of the index}

Once the Fredholm property is established, the index can change only if the family leaves the set of Fredholm operators. Compact perturbations and continuous deformations preserve this set; applied to the pseudodifferential calculus, they show that lower-order terms do not change the index.

\begin{theorem}[Compact perturbations and continuous families]
\label{pseudo:teo-estabilidad-indice}
\index{Fredholm index@Fredholm index!stability}
Let $X,Y$ be Hilbert spaces.
\begin{enumerate}[label=(\roman*)]
\item If $T\in\mathcal L(X,Y)$ is Fredholm and $K\in\mathcal K(X,Y)$, then $T+K$ is Fredholm and
\[
\operatorname{ind}(T+K)=\operatorname{ind}(T).
\]
\item The set of Fredholm operators is open in $\mathcal L(X,Y)$, and the index is locally constant.
\item If $\Lambda$ is connected and $\lambda\mapsto T_\lambda\in\mathcal L(X,Y)$ is a continuous family of Fredholm operators, then $\operatorname{ind}(T_\lambda)$ is constant on $\Lambda$.
\end{enumerate}
\end{theorem}

\begin{proof}
Since $T$ is Fredholm, there are orthogonal decompositions
\[
X=\ker T\oplus X_0,
\qquad
Y=\operatorname{Ran}T\oplus Y_0,
\]
with $Y_0$ finite-dimensional. The restriction $T\colon X_0\to\operatorname{Ran}T$ is a continuous isomorphism. Define $S\colon Y\to X$ as its inverse on $\operatorname{Ran}T$ and zero on $Y_0$. Then $ST-I_X$ and $TS-I_Y$ have finite rank. For $T+K$, the products with $S$ differ from the identity by compact operators, so Lemma~\ref{pseudo:lem-parametrix-compacta} proves that $T+K$ is Fredholm.

For operators sufficiently close to $T$, the components between $X_0$ and $\operatorname{Ran}T$ remain invertible by stability of invertibility. Eliminating this block by invertible operations reduces the nearby operator to a map between $\ker T$ and $Y_0$. Its index is
\[
\dim\ker T-\dim Y_0=\operatorname{ind}T,
\]
regardless of the rank of this finite block. This proves openness and local constancy simultaneously. Applied to the path $T+tK$, it proves equality of indices in (i). Finally, a locally constant integer-valued function is constant on a connected space, giving (iii).
\end{proof}

\begin{corollary}[Pseudodifferential stability]
\label{pseudo:cor-estabilidad-indice-pseudo}
Let $(M,\mathbf{g})$ be a smooth closed Riemannian manifold; let $\mathbf{E},\mathbf{F}\longrightarrow M$ be complex vector bundles of finite rank equipped with Hermitian bundle metrics and Hermitian connections; let $m,s\in\mathbb R$; and let $A\in\Psi^m_{\mathrm{cl}}(M;\mathbf{E},\mathbf{F})$ be elliptic.
\begin{enumerate}[label=(\alph*)]
\item If $K\colon H^{s+m}(M,\mathbf{E})\longrightarrow H^s(M,\mathbf{F})$ is compact, then $A+K$ is Fredholm with the same index as $A$.
\item In particular, this holds if $K$ is smoothing. It also holds for $K\in\Psi^{m-\varepsilon}(M;\mathbf{E},\mathbf{F})$ with $\varepsilon>0$, since
\[
H^{s+m}(M,\mathbf{E})\xrightarrow{\ K\ }H^{s+\varepsilon}(M,\mathbf{F})
\hookrightarrow H^s(M,\mathbf{F})
\]
and the second arrow is compact.
\item If $(A_\lambda)_{\lambda\in\Lambda}$ is an elliptic family of order $m$, continuous in the pseudodifferential topology of order $m$, and $\Lambda$ is connected, then $\operatorname{ind}A_\lambda$ is independent of $\lambda$.
\end{enumerate}
\end{corollary}

\begin{proof}
The first two parts follow from Proposition~\ref{pseudo:prop-regularizantes-compactos}, the preceding theorem, and Rellich. For (c), the mapping theorem shows that continuity in sufficiently many symbol seminorms implies continuity in
\[
\mathcal L\bigl(H^{s+m}(M,\mathbf{E}),H^s(M,\mathbf{F})\bigr).
\]
Each $A_\lambda$ is Fredholm by Theorem~\ref{pseudo:teo-fredholm-eliptico}; part (iii) of the preceding theorem then applies.
\end{proof}

The stability obtained here is entirely analytic: it uses parametrices, regularity, compactness, and norm continuity. On a closed manifold, these are exactly the arguments that turn ellipticity of the symbol into a well-defined stable integer.

\section{Complex powers, heat, and the analytic index}
\label{pseudo:sec-potencias-complejas-indice}

The parametrix of an elliptic operator inverts its symbol with respect to composition. To construct noninteger powers, we must uniformly invert an entire family: the resolvent $\lambda I-A$. The variable $\lambda$ does not have the same weight as the covariable $\xi$. If $A$ has order $m>0$, the natural scaling is
\[
(\xi,\lambda)\longmapsto(t\xi,t^m\lambda).
\]
This scaling governs the resolvent parametrix, complex powers, and their comparison with the spectral calculus.

Throughout this section, $(M,\mathbf{g})$ is a smooth closed Riemannian manifold of dimension $n$, and $\mathbf{E}\longrightarrow M$ is a complex vector bundle of finite rank equipped with a Hermitian bundle metric $\mathbf{h}_{\mathbf{E}}$ and a compatible connection $\nabla^{\mathbf{E}}$. These data remain fixed, and within this section we abbreviate
\[
\Gamma(\mathbf{E}):=\Gamma(M,\mathbf{E}),
\qquad
\Psi^\mu(M;\mathbf{E}):=\Psi^\mu(M;\mathbf{E},\mathbf{E}).
\]
Let \[
A\in\Psi^m_{\mathrm{cl}}(M;\mathbf{E}),\qquad m>0,
\] be a classical elliptic operator. Write $\displaystyle a\sim\displaystyle\sum_{j=0}^{\infty}a_{m-j}$ for its Fourier symbol in coordinates and local frames. If $A$ is differential, recall the conversion fixed at the beginning of the chapter:
\[
a_m(x,\xi)=i^m\boldsymbol{\sigma}_m(A)(x,\xi).
\]
The sectorial condition is formulated for $a_m$, since this is the symbol entering quantization. By $\boldsymbol{\sigma}(A)$ and $\rho(A)$ we mean the spectrum and resolvent of the closed realization on $L^2(M,\mathbf{E})$ with domain $H^m(M,\mathbf{E})$.

\begin{lemma}[Graph norm and elliptic realization]
\label{pseudo:lem-norma-grafica-realizacion-eliptica}
Let $(M,\mathbf{g})$ be a smooth closed Riemannian manifold; let $\mathbf{E}\longrightarrow M$ be a complex vector bundle of finite rank equipped with a Hermitian bundle metric $\mathbf{h}_{\mathbf{E}}$ and a compatible connection $\nabla^{\mathbf{E}}$; and let $A\in\Psi^m_{\mathrm{cl}}(M;\mathbf{E},\mathbf{E})$ be elliptic, with $m>0$. Then there are constants $c,C>0$ such that
\begin{equation}
 c\|\mathbf{u}\|_{H^m(M,\mathbf{E})}
 \leq \|\mathbf{u}\|_{L^2(M,\mathbf{E})}+\|A\mathbf{u}\|_{L^2(M,\mathbf{E})}
 \leq C\|\mathbf{u}\|_{H^m(M,\mathbf{E})},
 \qquad \mathbf{u}\in\Gamma(M,\mathbf{E}).
 \label{pseudo:eq-equivalencia-norma-grafica-A}
\end{equation}
The operator $A\colon H^m(M,\mathbf{E})\subset L^2(M,\mathbf{E})\longrightarrow L^2(M,\mathbf{E})$ is closed, $\Gamma(M,\mathbf{E})$ is a core, and its graph norm is equivalent to the $H^m(M,\mathbf{E})$ norm.
\end{lemma}

\begin{proof}
The mapping property $A:H^m(M,\mathbf{E})\to L^2(M,\mathbf{E})$ gives $\|\mathbf{u}\|_{L^2(M,\mathbf{E})}+\|A\mathbf{u}\|_{L^2(M,\mathbf{E})}\leq C\|\mathbf{u}\|_{H^m(M,\mathbf{E})}$. Let $B\in\Psi^{-m}_{\mathrm{cl}}(M;\mathbf{E})$ and $R\in\Psi^{-\infty}(M;\mathbf{E},\mathbf{E})$ satisfy $BA=I-R$. Then
\[
 \|\mathbf{u}\|_{H^m(M,\mathbf{E})}
 \leq\|BA\mathbf{u}\|_{H^m(M,\mathbf{E})}+\|R\mathbf{u}\|_{H^m(M,\mathbf{E})}
 \leq C_1\|A\mathbf{u}\|_{L^2(M,\mathbf{E})}+C_2\|\mathbf{u}\|_{L^2(M,\mathbf{E})},
\]
since $B:L^2(M,\mathbf{E})\to H^m(M,\mathbf{E})$ and $R:L^2(M,\mathbf{E})\to H^m(M,\mathbf{E})$ are continuous. This proves \eqref{pseudo:eq-equivalencia-norma-grafica-A} on $\Gamma(\mathbf{E})$ and, by density, on $H^m(M,\mathbf{E})$.

If $\mathbf{u}_j\in H^m(M,\mathbf{E})$, $\mathbf{u}_j\to \mathbf{u}$ in $L^2(M,\mathbf{E})$, and $A\mathbf{u}_j\to \mathbf{f}$ in $L^2(M,\mathbf{E})$, the first inequality applied to $\mathbf{u}_j-\mathbf{u}_k$ shows that $(\mathbf{u}_j)$ is Cauchy in $H^m(M,\mathbf{E})$. Its limit in $H^m(M,\mathbf{E})$ agrees with $\mathbf{u}$, and continuity $A:H^m(M,\mathbf{E})\to L^2(M,\mathbf{E})$ gives $A\mathbf{u}=\mathbf{f}$; hence the realization is closed. Finally, for $\mathbf{u}\in H^m(M,\mathbf{E})$ choose $\mathbf{v}_j\in\Gamma(\mathbf{E})$ with $\mathbf{v}_j\to \mathbf{u}$ in $H^m(M,\mathbf{E})$. The second inequality implies $\mathbf{v}_j\to \mathbf{u}$ in graph norm. This is exactly the assertion that $\Gamma(\mathbf{E})$ is a core.
\end{proof}

\begin{lemma}[Nonzero spectrum of a compact operator]
\label{pseudo:lem-espectro-compacto-general}
Let $K$ be a compact operator on a complex Hilbert space.
Every $\mu\in\sigma(K)\setminus\{0\}$ is an isolated eigenvalue.
For every $\varepsilon>0$ there are only finitely many distinct
eigenvalues with $|\mu|\geq\varepsilon$.
The generalized eigenspace
\[
 N_\mu:=\bigcup_{s\geq1}\ker(K-\mu I)^s
\]
equals $\ker(K-\mu I)^{s_0}$ for some $s_0$ and is finite
dimensional. Its dimension is the algebraic multiplicity of $\mu$.
\end{lemma}

\begin{proof}
For $\mu\neq0$, $K-\mu I$ is Fredholm of index zero, being a
compact perturbation of the isomorphism $-\mu I$.
If its kernel were zero, its cokernel would also be zero, and
the bounded inverse theorem would make it invertible.
Thus every nonzero spectral point is an eigenvalue,
and its kernel is finite dimensional.

Suppose that there are infinitely many distinct eigenvalues
$\mu_j$ with $|\mu_j|\geq\varepsilon>0$. Choose corresponding
eigenvectors, and let $Y_j$ be the span of the first $j$.
Eigenvectors of distinct eigenvalues are linearly independent,
and every $Y_j$ is invariant under $K$.
Choose a unit vector $e_j\in Y_j\cap Y_{j-1}^\perp$.
Since $(K-\mu_jI)Y_j\subset Y_{j-1}$, for $i<j$,
\[
 Ke_j-Ke_i=\mu_je_j+y_{ij},\qquad y_{ij}\in Y_{j-1}.
\]
Hence $\|Ke_j-Ke_i\|\geq|\mu_j|\geq\varepsilon$,
contradicting compactness. There are only finitely many
eigenvalues outside each disk about zero, and every nonzero
eigenvalue is therefore isolated.

Set $N_s:=\ker(K-\mu I)^s$.
Expanding the power gives
$(K-\mu I)^s=(-\mu)^sI+C_s$ with compact $C_s$, so $N_s$
is finite dimensional. If $N_s=N_{s+1}$, all subsequent
spaces coincide: for $x\in N_{s+2}$,
$(K-\mu I)x\in N_{s+1}=N_s$, hence $x\in N_{s+1}$;
continue by induction.
If no equality occurred, we could choose unit vectors
$e_s\in N_s\cap N_{s-1}^\perp$ for every $s$, with
$N_0=\{0\}$. Since $(K-\mu I)e_s\in N_{s-1}$ and
$Ke_i\in N_{s-1}$ for $i<s$, we would have
$\|Ke_s-Ke_i\|\geq|\mu|$.
Compactness excludes this possibility.
The chain stabilizes, and its union is finite dimensional.
\end{proof}

\begin{proposition}[Compact resolvent and discrete spectrum]
\label{pseudo:prop-resolvente-compacto-espectro-discreto}
Let $(M,\mathbf{g})$ be a smooth closed Riemannian manifold; let $\mathbf{E}\longrightarrow M$ be a complex vector bundle of finite rank equipped with a Hermitian bundle metric $\mathbf{h}_{\mathbf{E}}$ and a compatible connection $\nabla^{\mathbf{E}}$; and let $A\in\Psi^m_{\mathrm{cl}}(M;\mathbf{E},\mathbf{E})$ be elliptic, with $m>0$, realized on $L^2(M,\mathbf{E})$ with domain $H^m(M,\mathbf{E})$. If $\rho(A)\neq\varnothing$, every resolvent of this realization is compact: for each $\lambda_0\in\rho(A)$,
\[
(\lambda_0I-A)^{-1}:L^2(M,\mathbf{E})\longrightarrow L^2(M,\mathbf{E})
\]
is compact. Under this hypothesis, $\sigma(A)$ consists of isolated eigenvalues, each of finite algebraic multiplicity, and has no finite accumulation points. If $A$ is self-adjoint and bounded below, then $\rho(A)\neq\varnothing$ and there is an orthonormal basis $(\boldsymbol{\phi}_j)$ of $L^2(M,\mathbf{E})$ consisting of smooth eigensections. Repeating eigenvalues according to multiplicity, we have
\[
 A\boldsymbol{\phi}_j=\lambda_j\boldsymbol{\phi}_j.
\]
If \(L^2(M,\mathbf E)\) is infinite dimensional, the basis can be
enumerated with \(\lambda_j\to+\infty\). If that space has finite
dimension, the basis and the list of eigenvalues are finite;
in dimension zero both are empty. Spectral sums in this section
are interpreted as finite sums in that case.

\end{proposition}

\begin{proof}
A parametrix $B\in\Psi^{-m}_{\mathrm{cl}}(M;\mathbf{E})$ and a remainder $R\in\Psi^{-\infty}(M;\mathbf{E},\mathbf{E})$ give $\mathbf{u}=BA\mathbf{u}+R\mathbf{u}$ for $\mathbf{u}\in\Gamma(M,\mathbf{E})$. The mapping properties and extension by density yield the elliptic estimate
\begin{equation}
\|\mathbf{u}\|_{H^m(M,\mathbf{E})}
\leq C\bigl(\|A\mathbf{u}\|_{L^2(M,\mathbf{E})}+\|\mathbf{u}\|_{L^2(M,\mathbf{E})}\bigr),
\qquad \mathbf{u}\in H^m(M,\mathbf{E}).
\label{pseudo:eq-estimacion-eliptica-resolvente-compacto}
\end{equation}
If $\mathbf{u}=(\lambda_0I-A)^{-1}f$, continuity of the resolvent on $L^2(M,\mathbf{E})$ and $A\mathbf{u}=\lambda_0\mathbf{u}-f$ show, through \eqref{pseudo:eq-estimacion-eliptica-resolvente-compacto}, that $(\lambda_0I-A)^{-1}$ is continuous from $L^2(M,\mathbf{E})$ to $H^m(M,\mathbf{E})$. The inclusion $H^m(M,\mathbf{E})\hookrightarrow L^2(M,\mathbf{E})$ is compact by Rellich, so the resolvent viewed on $L^2(M,\mathbf{E})$ is compact.

The spectral correspondence
\[
\lambda\in\sigma(A)
\quad\Longleftrightarrow\quad
(\lambda_0-\lambda)^{-1}
 \in\sigma\bigl((\lambda_0I-A)^{-1}\bigr),
\qquad \lambda\neq\lambda_0,
\]
follows from the factorization as maps from $\mathcal D(A)$ to $L^2$,
\[
 \lambda I-A
 =\bigl(I+(\lambda-\lambda_0)K\bigr)(\lambda_0I-A),
 \qquad K:=(\lambda_0I-A)^{-1}.
\]
Lemma~\ref{pseudo:lem-espectro-compacto-general} gives isolation
of nonzero spectral points of $K$, and hence of the corresponding
eigenvalues of $A$.
To identify multiplicity, let
$\mu=(\lambda_0-\lambda)^{-1}$ and let $N_\mu$ be the generalized
eigenspace of $K$. Then $K|_{N_\mu}=\mu I+N$ with nilpotent $N$.
This restriction is invertible, and
$N_\mu=K(N_\mu)\subset\mathcal D(A)$.
On $N_\mu$,
\[
 A=\lambda_0I-(K|_{N_\mu})^{-1},\qquad
 K-\mu I=\mu K(A-\lambda I).
\]
The factors $\mu K$ are invertible and commute with $A-\lambda I$,
so their powers have the same kernels.
Conversely, for $u\in\ker(A-\lambda I)^s$, iterating the same
identity gives $(K-\mu I)^su=0$.
Thus the generalized eigenspaces of $A$ and $K$ coincide
at corresponding values, and their dimension is finite.

In the self-adjoint bounded-below case, choose a real $\lambda_0$ below the spectrum. The resolvent is then compact and self-adjoint, and its spectral theorem provides an orthonormal basis of eigenvectors. If the space is infinite dimensional, the nonzero resolvent
eigenvalues tend to zero; the preceding correspondence and the lower
bound then force \(\lambda_j\to+\infty\). In finite dimension the
spectral theorem gives a finite basis and no such divergence is asserted.

Finally, $A\boldsymbol{\phi}_j=\lambda_j\boldsymbol{\phi}_j$ and elliptic regularity imply $\boldsymbol{\phi}_j\in\Gamma(\mathbf{E})$.
\end{proof}

\subsection{Spectral cut and a minimal parameter-dependent calculus}

Constructing $A^z$ requires choosing a branch of the logarithm and controlling the resolvent along a sector avoiding both the principal spectrum and the spectrum of the closed realization. We first fix the cut, then construct the uniform parametrix of $\lambda I-A$.

For $\theta\in\mathbb R$, set
\[
L_\theta:=\{re^{i\theta}\mid r>0\}.
\]
On $\mathbb C\setminus(L_\theta\cup\{0\})$, fix the branch \[
\log_\theta\lambda
=\log|\lambda|+i\arg_\theta\lambda,
\qquad
\theta-2\pi<\arg_\theta\lambda<\theta,
\] and write $\lambda^z:=e^{z\log_\theta\lambda}$.

\begin{definition}[Agmon angle and spectral cut]
\label{pseudo:def-angulo-agmon-corte}
\index{Agmon angle@Agmon angle}
\index{spectral cut}
Let $(M,\mathbf{g})$ be a smooth closed Riemannian manifold, let $(\mathbf{E},\mathbf{h}_{\mathbf{E}})\longrightarrow M$ be a complex vector bundle of finite rank equipped with a Hermitian bundle metric, and let $A\in\Psi^m_{\mathrm{cl}}(M;\mathbf{E},\mathbf{E})$ be elliptic, with $m>0$ and principal Fourier symbol $a_m$. Let $\theta\in\mathbb R$. We say that $\theta$ is an \textbf{Agmon angle for the principal symbol} of $A$ if
\[
\operatorname{spec}a_m(x,\xi)\cap L_\theta=\varnothing
\]
for every $(x,\xi)\in T^*M\setminus0_M$. It is an \textbf{admissible spectral cut for $A$} if, in addition,
\[
\sigma(A)\cap L_\theta=\varnothing.
\]
To define all powers, positive and negative, we also assume $0\in\rho(A)$, that is, invertibility of $A$. If $0$ is an eigenvalue, reduced powers require separating the generalized eigenspace, as specified next.
\end{definition}

If $0\in\sigma(A)$ and $A$ is not necessarily self-adjoint, the correct separation uses more than $\ker A$ alone. Existence of the admissible spectral cut implies $\rho(A)\neq\varnothing$; Proposition~\ref{pseudo:prop-resolvente-compacto-espectro-discreto} then shows that $0$ is isolated. Let \begin{equation}
\Pi_0
:=
\frac{1}{2\pi i}\int_{|\lambda|=\varepsilon}
(\lambda I-A)^{-1}\,d\lambda
\label{pseudo:eq-proyeccion-riesz-cero}
\end{equation} be its Riesz projection, with the circle positively oriented and containing no other spectral point. Its range is the finite-dimensional generalized eigenspace of $0$, and $L^2(M,\mathbf{E})=\operatorname{ran}\Pi_0\oplus\ker\Pi_0$ is a $A$-invariant decomposition. The restriction of $A$ to $\ker\Pi_0$ is invertible; by \textbf{reduced power} we mean exclusively the power of this restriction, extended by zero on $\operatorname{ran}\Pi_0$. This convention asserts no group law on all of $L^2(M,\mathbf{E})$. If $A$ is self-adjoint, $\Pi_0$ is the orthogonal projection onto $\ker A$, and no nilpotent part occurs.

In a cotangent chart, with the Euclidean norm fixed by this chapter's local convention, the quantity \[
\rho_m(\xi,\lambda):=\|\xi\|_{\mathbb R^n}+|\lambda|^{\frac{1}{m}}
\] is homogeneous of degree one under the weighted scaling fixed above. For smooth excisions, use the equivalent weight
\begin{equation}
\widehat\rho_m(\xi,\lambda)
:=
\bigl(\langle\xi\rangle^{2m}+|\lambda|^2\bigr)^{\frac{1}{2m}}.
\label{pseudo:eq-peso-suave-parametro}
\end{equation}
Unlike $\rho_m$, this weight is smooth also on the axes; the two are comparable in the region of high frequencies and large parameters.

\begin{definition}[Parameter-dependent symbol and ellipticity]
\label{pseudo:def-simbolo-elipticidad-parametro}
\index{parameter-dependent symbol@parameter-dependent symbol}
\index{ellipticity!with parameter}
Let $m>0$, let $U\subseteq\mathbb R^n$ be open, let $r\in\mathbb N$, and let $\Lambda\subseteq\mathbb C$ be a closed sector. Consider families with values in $\operatorname{End}(\mathbb C^r)$ and set
\[
\Omega_m(\xi,\lambda):=|\lambda|+\langle\xi\rangle^m.
\]
If \[
c\in C^\infty\!\left(
U\times\mathbb R^n\times\Lambda;
\operatorname{End}(\mathbb C^r)
\right)
\] is smooth up to the faces of the sector in the real variables $(x,\xi,\operatorname{Re}\lambda,\operatorname{Im}\lambda)$, for $K\Subset U$, $j\in\mathbb N_0$, $\alpha,\beta\in\mathbb N_0^n$, and $k,\ell\in\mathbb N_0$ define the following seminorms, where $\|\cdot\|_{\operatorname{op},\mathbb C^r}$ is the operator norm induced by the standard Hermitian norm:
\begin{align}
 \mathfrak p^{\mathrm R,j}_{K,\alpha,\beta,k,\ell}(c)
 &:={}
 \sup_{\substack{x\in K,\,\xi\in\mathbb R^n\\\lambda\in\Lambda}}
 \langle\xi\rangle^{j+|\beta|}
 \Omega_m(\xi,\lambda)^{1+k+\ell}
 \bigl\|D_x^\alpha D_\xi^\beta
 \partial_\lambda^k\partial_{\bar\lambda}^{\ell}c(x,\xi,\lambda)\bigr\|_{\operatorname{op},\mathbb C^r},
 \label{pseudo:eq-seminorma-debil-resolvente}\\
 \mathfrak p^{\mathrm D,j}_{K,\alpha,\beta,k,\ell}(c)
 &:={}
 \sup_{\substack{x\in K,\,\xi\in\mathbb R^n\\\lambda\in\Lambda}}
 \langle\xi\rangle^{j+|\beta|}
 \Omega_m(\xi,\lambda)^{k+\ell}
 \bigl\|D_x^\alpha D_\xi^\beta
 \partial_\lambda^k\partial_{\bar\lambda}^{\ell}c(x,\xi,\lambda)\bigr\|_{\operatorname{op},\mathbb C^r}.
 \label{pseudo:eq-seminorma-debil-defecto}
\end{align}
Here $D_x^\alpha=\partial_x^\alpha$ and $D_\xi^\beta=\partial_\xi^\beta$, in accordance with this chapter's convention, while $\partial_\lambda=\displaystyle\frac{1}{2}
(\partial_{\operatorname{Re}\lambda}
-i\partial_{\operatorname{Im}\lambda})$ and $\partial_{\bar\lambda}=\displaystyle\frac{1}{2}
(\partial_{\operatorname{Re}\lambda}
+i\partial_{\operatorname{Im}\lambda})$ are the Wirtinger derivatives.

A family $q^{(j)}$, with $j\geq0$, has \textbf{resolvent type $j$} if all seminorms \eqref{pseudo:eq-seminorma-debil-resolvente} are finite. It has \textbf{defect type $j$} if all seminorms \eqref{pseudo:eq-seminorma-debil-defecto} are finite. In particular, for the first type,
\[
 \|D_x^\alpha D_\xi^\beta\partial_\lambda^k
 \partial_{\bar\lambda}^\ell q^{(j)}(x,\xi,\lambda)\|_{\operatorname{op},\mathbb C^r}
 \leq C_{K,\alpha,\beta,k,\ell,j}
 \langle\xi\rangle^{-j-|\beta|}
 \Omega_m(\xi,\lambda)^{-1-k-\ell},
\]
and for the second the last factor is replaced by $\Omega_m^{-k-\ell}$. Families of resolvent type $j$ and defect type $j$ behave as orders $-m-j$ and $-j$, respectively, in the cone $|\lambda|\leq C\langle\xi\rangle^m$, while retaining the correct estimate when $|\lambda|$ dominates.

A family is \textbf{holomorphic in the parameter} if $\partial_{\bar\lambda}c=0$ in the interior of $\Lambda$. It is \textbf{simultaneously smoothing in $\xi$ and $\lambda$} if, for all $P,Q\in\mathbb N_0$, every compact set $K$, and all $\alpha,\beta,k,\ell$, there exists a constant $C_{K,P,Q,\alpha,\beta,k,\ell}>0$ such that
\begin{equation}
 \bigl\|D_x^\alpha D_\xi^\beta\partial_\lambda^k
 \partial_{\bar\lambda}^{\ell}c(x,\xi,\lambda)\bigr\|_{\operatorname{op},\mathbb C^r}
 \leq C_{K,P,Q,\alpha,\beta,k,\ell}
 \langle\xi\rangle^{-P}(1+|\lambda|)^{-Q}.
 \label{pseudo:eq-regularizacion-simultanea-parametro}
\end{equation}
We abbreviate this condition by saying that $c$ is smoothing with parameter. Membership in every resolvent type $j$ imposes only arbitrary decay in $\xi$ and the resolvent factor in $\lambda$; condition \eqref{pseudo:eq-regularizacion-simultanea-parametro} additionally requires every power of $1+|\lambda|$ and does not follow from the former condition.

This is the weak version of the parameter-dependent calculus. The strong bound condensed into a single power of $\rho_m$ fails for a general classical pseudodifferential operator when $|\lambda|$ dominates $\langle\xi\rangle^m$; the preceding separate bounds remain valid. At high frequencies, a family of resolvent type $0$ is \textbf{classical with parameter} if its components are homogeneous under $(\xi,\lambda)\mapsto(t\xi,t^m\lambda)$ and the remainder after $N$ components has resolvent type $N$ for $\|\xi\|\geq1$. The complementary bounded-frequency part is required only to be smoothing in $\xi$, retaining the decay factor $\Omega_m^{-1-k-\ell}$. If $c_j$ has resolvent type $j$, write
\begin{equation}
 c\sim\sum_{j=0}^{\infty}c_j
 \quad\Longleftrightarrow\quad
 c-\sum_{j=0}^{N-1}c_j\ \text{has resolvent type $N$ for every }
 N\in\mathbb N_0.
 \label{pseudo:eq-def-suma-asintotica-debil}
\end{equation}
The same convention, replacing \emph{resolvent} by \emph{defect}, is used for defect-type families.

The homogeneous principal symbol $a_m$ is \textbf{elliptic with parameter} in $\Lambda$ if $\lambda I-a_m(x,\xi)$ is invertible for $\xi\neq0$ and $\lambda\in\Lambda$ and, for every $K\Subset U$, there exists $C_K>0$ such that its inverse satisfies
\[
 \|\bigl(\lambda I-a_m(x,\xi)\bigr)^{-1}\|
 \leq C_K\bigl(|\lambda|+\|\xi\|^m\bigr)^{-1},
 \qquad x\in K.
\]
This condition is stronger than ellipticity of $A$ without a parameter.
\end{definition}

The smooth weight in \eqref{pseudo:eq-peso-suave-parametro} satisfies
\begin{equation}
 2^{-\frac{1}{2}}\Omega_m(\xi,\lambda)
 \leq \widehat\rho_m(\xi,\lambda)^m
 \leq \Omega_m(\xi,\lambda).
 \label{pseudo:eq-comparacion-pesos-parametro}
\end{equation}
Indeed, for $a,b\geq0$ we have $\displaystyle\frac{a+b}{\sqrt2}\leq(a^2+b^2)^{\displaystyle\frac{1}{2}}\leq a+b$; it suffices to take $a=\langle\xi\rangle^m$ and $b=|\lambda|$. If $\chi\in C^\infty(\mathbb R)$ is constant outside a compact subset of $(0,\infty)$, fix $0<c_\chi<C_\chi$ such that $\operatorname{supp}\chi'\subseteq[c_\chi,C_\chi]$. The chain rule gives, when $|\gamma|+r+s\geq1$,
\begin{equation}
 \bigl|D_\xi^\gamma\partial_\lambda^r
 \partial_{\bar\lambda}^{s}\chi(\varepsilon\widehat\rho_m)\bigr|
 \leq C_{\gamma,r,s}
 \langle\xi\rangle^{-|\gamma|}\Omega_m^{-r-s}
 \mathbf 1_{\{c_\chi\varepsilon^{-1}
 \leq\widehat\rho_m\leq C_\chi\varepsilon^{-1}\}}.
 \label{pseudo:eq-derivadas-corte-rho-parametro}
\end{equation}
To justify every derivative order, first differentiate $\widehat\rho_m^{2m}=\langle\xi\rangle^{2m}+|\lambda|^2$. A derivative in $\xi$ of the resulting quotient is bounded by $C\langle\xi\rangle^{-1}$, and a Wirtinger derivative by $C\widehat\rho_m^{-m}$. Repeated use of the Leibniz rule gives terms containing $|\gamma|$ factors of the first type and $r+s$ of the second, possibly grouped into higher-order derivatives; their total exponents remain the same. Comparison \eqref{pseudo:eq-comparacion-pesos-parametro} proves \eqref{pseudo:eq-derivadas-corte-rho-parametro}.

\begin{lemma}[Asymptotic summation in weak parameter-dependent classes]
\label{pseudo:lem-suma-asintotica-parametro-debil}
Fix $m>0$, an open set $U\subseteq\mathbb R^n$, a rank $r\in\mathbb N$, and a closed sector $\Lambda\subseteq\mathbb C$ as in Definition~\ref{pseudo:def-simbolo-elipticidad-parametro}. Let $(c_j)_{j\geq0}$ be a sequence of $\operatorname{End}(\mathbb C^r)$-valued families such that $c_j$ has resolvent type $j$. There is a family $c$ of resolvent type $0$ satisfying
$\displaystyle c\sim\displaystyle\sum_{j=0}^{\infty}c_j$
in the sense of \eqref{pseudo:eq-def-suma-asintotica-debil}. The assertion also holds for defect-type families. If all $c_j$ are holomorphic in $\lambda$, the sum can be chosen holomorphic. Any two sums differ by a family belonging to every resolvent type $N$; the difference is smoothing with parameter precisely when it additionally satisfies the powers of $1+|\lambda|$ in \eqref{pseudo:eq-regularizacion-simultanea-parametro}. Finally, suppose that, for each $P,Q$ and each finite family of derivatives, there exists $J$ such that the terms $c_j$, $j\geq J$, satisfy an inequality of the form \eqref{pseudo:eq-regularizacion-simultanea-parametro} with constants that may depend on $j$. Then one may use the joint cutoff $\chi(\varepsilon_j\widehat\rho_m)$ and choose $\varepsilon_j$ with the following quantified property. For each simultaneous seminorm $\mathfrak s_{K,P,Q,\alpha,\beta,k,\ell}$, there exists $J$ such that, for every $N\geq J$,
\[
 \mathfrak s_{K,P,Q,\alpha,\beta,k,\ell}
 \left(\sum_{j=N}^{\infty}
 \chi(\varepsilon_j\widehat\rho_m)c_j\right)
 \leq \sum_{j=N}^{\infty}2^{-j}.
\]
This second auxiliary sum is not generally holomorphic in $\lambda$.
\end{lemma}

\begin{proof}
Choose $\chi\in C^\infty(\mathbb R)$ with $0\leq\chi\leq1$, $\chi(t)=0$ for $t\leq1$, and $\chi(t)=1$ for $t\geq2$. Exhaust each chart by compact sets $K_1\Subset K_2\Subset\cdots$ and enumerate in a sequence $(\mathfrak p_\nu)_{\nu\geq1}$ the seminorms \eqref{pseudo:eq-seminorma-debil-resolvente} with $K=K_d$ and all indices $d,\alpha,\beta,k,\ell$. For $0<\varepsilon\leq1$, set
\[
 \chi_{j,\varepsilon}(\xi):=\chi(\varepsilon\langle\xi\rangle).
\]
Successive applications of the chain rule give
\begin{equation}
 |D_\xi^\gamma\chi_{j,\varepsilon}(\xi)|
 \leq C_\gamma\langle\xi\rangle^{-|\gamma|}
 \mathbf 1_{\{\langle\xi\rangle\geq\varepsilon^{-1}\}}.
 \label{pseudo:eq-derivadas-corte-suma-debil}
\end{equation}
For $\gamma=0$, the indicator merely expresses the cutoff support. Where a derivative of the cutoff is nonzero, $\varepsilon\langle\xi\rangle\in[1,2]$; each derivative of $\langle\xi\rangle$ contributes a power bounded by that on the right-hand side.

The Leibniz rule, \eqref{pseudo:eq-derivadas-corte-suma-debil}, and the defining inequality for $c_j$ show that each term of $D_\xi^\beta(\chi_{j,\varepsilon}c_j)$ is bounded by
\[
 C\langle\xi\rangle^{-j-|\beta|}
 \Omega_m^{-1-k-\ell}
 \mathbf 1_{\{\langle\xi\rangle\geq\varepsilon^{-1}\}}.
\]
Thus, for each seminorm $\mathfrak p_\nu$ of resolvent type $0$, there exists $C_{j,\nu}>0$ such that
\begin{equation}
 \mathfrak p_\nu(\chi_{j,\varepsilon}c_j)
 \leq C_{j,\nu}\varepsilon^j.
 \label{pseudo:eq-control-corte-seminorma-parametro}
\end{equation}
The derivatives $\partial_\lambda^k\partial_{\bar\lambda}^{\ell}$ fall only on $c_j$, since this cutoff depends on $\xi$ and not on $\lambda$; hence the exponent $-1-k-\ell$ of $\Omega_m$ is preserved exactly.

Set $\varepsilon_0=1$ and recursively choose $\varepsilon_j\in(0,2^{-j}]$ for $j\geq1$. At step $j$, require, for the first $j$ seminorms and every $N=0,1,\ldots,j-1$, that the resolvent-type $N$ seminorm of the cutoff term not exceed $2^{-j}$; for $N=0$ this is equivalent to requiring
\begin{equation}
 C_{j,\nu}\varepsilon_j^j\leq2^{-j},
 \qquad 1\leq\nu\leq j.
 \label{pseudo:eq-eleccion-epsilon-suma-parametro}
\end{equation}
There are only finitely many conditions at each step, all containing a positive power $\varepsilon_j^{j-N}$, so the choice is possible. Define
\begin{equation}
 c:=\sum_{j=0}^{\infty}
 \chi(\varepsilon_j\langle\xi\rangle)c_j.
 \label{pseudo:eq-suma-borel-parametro}
\end{equation}
The series is locally finite: for fixed $\xi$, $\varepsilon_j\langle\xi\rangle\leq2^{-j}\langle\xi\rangle<1$ holds from some index onward. For a fixed seminorm $\mathfrak p_\nu$, its tail from $j\geq\displaystyle\max\{\nu,1\}$ is dominated by $\displaystyle \sum_{j=0}^{\infty}2^{-j}$; hence it converges in every seminorm of resolvent type $0$.

Fix $N$. For $j>N$, the preceding calculation, normalized with $\langle\xi\rangle^{N+|\beta|}$, contains $\langle\xi\rangle^{N-j}\leq\varepsilon_j^{j-N}$ on the cutoff support. The conditions imposed when choosing $\varepsilon_j$ prove convergence of $\displaystyle \sum_{j=N}^{\infty}\chi(\varepsilon_j\langle\xi\rangle)c_j$ in every seminorm of resolvent type $N$. Moreover,
\begin{equation}
 c-\sum_{j=0}^{N-1}c_j
 =\sum_{j=N}^{\infty}\chi(\varepsilon_j\langle\xi\rangle)c_j
 -\sum_{j=0}^{N-1}
 \bigl(1-\chi(\varepsilon_j\langle\xi\rangle)\bigr)c_j.
 \label{pseudo:eq-resto-suma-asintotica-parametro}
\end{equation}
The second sum is finite, and each summand has $\xi$ in a compact set, so it belongs to every resolvent type $N$. This proves the remainder estimate for each $N$.

For defect type, repeat the same inequalities with $\Omega_m^{-k-\ell}$; no step changes that power. If the $c_j$ are holomorphic, the locally finite sum \eqref{pseudo:eq-suma-borel-parametro} is also holomorphic. Subtracting two sums and using \eqref{pseudo:eq-resto-suma-asintotica-parametro} for every $N$ shows that their difference belongs to all the indicated types. If it additionally has the factor $(1+|\lambda|)^{-Q}$ for each $Q$, one obtains exactly \eqref{pseudo:eq-regularizacion-simultanea-parametro}; this is the quantified formulation of simultaneous residual decay. To prove the last assertion, now use $\chi(\varepsilon_j\widehat\rho_m)$ and also enumerate the seminorms
\[
 \mathfrak s_{K,P,Q,\alpha,\beta,k,\ell}(d)
 :=\sup_{\substack{x\in K,\,\xi\in\mathbb R^n\\\lambda\in\Lambda}}
 \langle\xi\rangle^P(1+|\lambda|)^Q
 \|D_x^\alpha D_\xi^\beta\partial_\lambda^k
 \partial_{\bar\lambda}^{\ell}d\|.
\]
Fix one of them and apply the hypothesis to exponents $P+L$ and $Q+L$. The Leibniz rule and \eqref{pseudo:eq-derivadas-corte-rho-parametro} give, on the cutoff support, the bound
\[
 \mathfrak s_{K,P,Q,\alpha,\beta,k,\ell}
 \bigl(\chi(\varepsilon_j\widehat\rho_m)c_j\bigr)
 \leq C_{j,L}\sup_{\substack{x\in K,\ \xi\in\mathbb R^n,\ \lambda\in\Lambda\\\widehat\rho_m(\xi,\lambda)\geq\varepsilon_j^{-1}}}
 \langle\xi\rangle^{-L}(1+|\lambda|)^{-L}.
\]
By \eqref{pseudo:eq-comparacion-pesos-parametro}, every point in the last supremum satisfies $\langle\xi\rangle\geq c\widehat\rho_m$ or $1+|\lambda|\geq c\widehat\rho_m^m$. Thus the right-hand side does not exceed $C_{j,L}\varepsilon_j^{\min\{1,m\}L}$. For a weak seminorm of resolvent or defect type, first absorb the positive powers of $\langle\xi\rangle$ and $\Omega_m$ by choosing $P+L$ and $Q+L$ still larger, obtaining the same positive power of $\varepsilon_j$.

Recursively decrease $\varepsilon_j$ so that the first $j$ applicable seminorms do not exceed $2^{-j}$. For a fixed seminorm, increase the index $J$ in the hypothesis if necessary until that seminorm is among the first $j$ for every $j\geq J$. The condition is then imposed on each summand of index $j\geq J$, and for every $N\geq J$ its tail is dominated by $\displaystyle \sum_{j=N}^{\infty}2^{-j}$. Moreover, $1-\chi(\varepsilon_j\widehat\rho_m)$ has jointly compact support in $(\xi,\lambda)$ for each fixed $j$. Thus the finite corrections are smoothing with parameter, and the tails satisfy \eqref{pseudo:eq-regularizacion-simultanea-parametro}. This proves the assertion without identifying simultaneous smoothing with mere membership in every resolvent type.
\end{proof}

The cutoff in this proof depends only on $\xi$ and preserves holomorphy. In contrast, the global cutoffs $\chi(\varepsilon\widehat\rho_m)$ needed to excise the origin $(\xi,\lambda)=(0,0)$ jointly contain $|\lambda|^2$ and, unless constant, are not holomorphic. We always distinguish the auxiliary symbol family containing these cutoffs, the holomorphic homogeneous components, and the holomorphic exact resolvent.

The first condition is not a statement about the global spectrum of $A$. It is a uniform condition on a finite-dimensional family of matrices. Homogeneity of $a_m$ and compactness of the cosphere bundle allow the ray to be widened: there exists $\varepsilon_{\mathrm p}>0$ such that
\[
\lambda I-a_m(x,\xi)
\quad\text{is invertible if}\quad
\xi\neq0,
\quad
\lambda\in\Lambda^{\mathrm p}_{\theta,\varepsilon_{\mathrm p}},
\]
where \[
\Lambda^{\mathrm p}_{\theta,\varepsilon_{\mathrm p}}
:=
\{re^{i\varphi}\mid r>0,\ |\varphi-\theta|
\leq\varepsilon_{\mathrm p}\}
\] is understood modulo $2\pi$. We do not yet assert that this principal sector is contained in $\rho(A)$. That conclusion is obtained first for large $|\lambda|$ using the parametrix and Neumann series. Only afterward is the sector narrowed to avoid the remaining compact spectral set; this logical order avoids using resolvent existence to construct the resolvent. Until then, all subsequent parameter seminorms and asymptotic sums are taken with $\Lambda=\Lambda^{\mathrm p}_{\theta,\varepsilon_{\mathrm p}}$.

Inversion on the normalized compact set followed by dilation gives
\[
\bigl\|(\lambda I-a_m(x,\xi))^{-1}\bigr\|_{\mathcal L((\mathbf{E}_x,\mathbf{h}_{\mathbf{E}}(x)))}
\leq C\bigl(|\lambda|+|\xi|_{\mathbf{g}}^m\bigr)^{-1}.
\]
Differentiating the identity $(\lambda I-a_m)(\lambda I-a_m)^{-1}=I$ gives, for the homogeneous terms introduced below and, in a local unitary frame, for $\|\xi\|_{\mathbb R^n}\geq1$, weak estimates of the form
\[
\left\|
D_x^\alpha D_\xi^\beta\partial_\lambda^kq_{-m-j}(x,\xi,\lambda)
\right\|_{\operatorname{op},\mathbb C^r}
\leq
C_{\alpha,\beta,k}
\langle\xi\rangle^{-j-|\beta|}
\Omega_m(\xi,\lambda)^{-1-k}.
\]
These bounds, together with holomorphy in $\lambda$, are the portion of the parameter-dependent calculus we use.

\begin{proposition}[Resolvent parametrix]
\label{pseudo:prop-parametrix-resolvente}
Let $(M,\mathbf{g})$ be a smooth closed Riemannian manifold; let $(\mathbf{E},\mathbf{h}_{\mathbf{E}})\longrightarrow M$ be a complex vector bundle of finite rank equipped with a Hermitian bundle metric $\mathbf{h}_{\mathbf{E}}$ and a compatible connection $\nabla^{\mathbf{E}}$; and let $A\in\Psi^m_{\mathrm{cl}}(M;\mathbf{E},\mathbf{E})$ be elliptic, with $m>0$ and local full symbol $a$. Suppose $\theta$ is an admissible spectral cut for $A$ and $0\in\rho(A)$. There exists a global auxiliary symbol family $\widetilde q$ of resolvent type $0$, not necessarily holomorphic, which is a two-sided inverse of $\lambda I-a$ modulo families smoothing with parameter in the principal sector $\Lambda^{\mathrm p}_{\theta,\varepsilon_{\mathrm p}}$. At high frequencies its classical expansion is
\[
\widetilde q(x,\xi,\lambda)
\underset{\|\xi\|_{\mathbb R^n}\to\infty}{\sim}
\sum_{j=0}^\infty q_{-m-j}(x,\xi,\lambda),
\qquad
\lambda\in\Lambda^{\mathrm p}_{\theta,\varepsilon_{\mathrm p}},
\]
where the components are considered only for $\xi\neq0$, are holomorphic in $\lambda$, satisfy
\[
q_{-m-j}(x,t\xi,t^m\lambda)
=t^{-m-j}q_{-m-j}(x,\xi,\lambda),
\qquad t>0,
\]
and are determined recursively as follows. If \[
p_m:=\lambda I-a_m,
\qquad
p_{m-k}:=-a_{m-k}\quad(k\geq1),
\], then
\[
q_{-m}=p_m^{-1}
\]
and, for $j\geq1$,
\begin{equation}
q_{-m-j}
=
-q_{-m}
\sum_{\substack{k+\ell+|\alpha|=j\\0\leq\ell<j}}
\frac{i^{-|\alpha|}}{\alpha!}
(D_\xi^\alpha p_{m-k})
(D_x^\alpha q_{-m-\ell}).
\label{pseudo:eq-recurrencia-parametrix-resolvente}
\end{equation}
The same component also satisfies the right recurrence
\begin{equation}
 q_{-m-j}
 =-
 \sum_{\substack{k+\ell+|\alpha|=j\\0\leq\ell<j}}
 \frac{i^{-|\alpha|}}{\alpha!}
 (D_\xi^\alpha q_{-m-\ell})
 (D_x^\alpha p_{m-k})q_{-m}.
 \label{pseudo:eq-recurrencia-derecha-parametrix-resolvente}
\end{equation}
For every compact set $K$, all $\alpha,\beta\in\mathbb N_0^n$, and $k\in\mathbb N_0$, there exists $C_{K,j,\alpha,\beta,k}>0$ such that, if $\|\xi\|\geq1$,
\begin{equation}
 \bigl\|D_x^\alpha D_\xi^\beta\partial_\lambda^k
 q_{-m-j}(x,\xi,\lambda)\bigr\|
 \leq C_{K,j,\alpha,\beta,k}
 \langle\xi\rangle^{-j-|\beta|}
 \Omega_m(\xi,\lambda)^{-1-k},
 \qquad
 \partial_{\bar\lambda}q_{-m-j}=0.
 \label{pseudo:eq-cotas-componentes-resolvente}
\end{equation}
There exist $0<\varepsilon\leq\varepsilon_{\mathrm p}$ and $R_0>0$ such that
\[
 \Lambda_{\theta,\varepsilon}
 :=\{re^{i\varphi}\mid r>0,\ |\varphi-\theta|\leq\varepsilon\}
 \subseteq\rho(A).
\]
The full symbol of the resolvent in this sector, $r(x,\xi,\lambda)$, taken in a fixed local quantization, is holomorphic in $\lambda$; for $|\lambda|\geq R_0$ it has the same high-frequency expansion
\[
r(x,\xi,\lambda)
\underset{\|\xi\|\to\infty}{\sim}
\sum_{j=0}^\infty q_{-m-j}(x,\xi,\lambda).
\]
For every $s\in\mathbb R$, there is a constant $C_s>0$ such that, for $\lambda\in\Lambda_{\theta,\varepsilon}$ with $|\lambda|\geq R_0$, the resolvent satisfies
\begin{equation}
\|(\lambda I-A)^{-1}\mathbf{u}\|_{H^{s+m}(M,\mathbf{E})}
+|\lambda|\|(\lambda I-A)^{-1}\mathbf{u}\|_{H^s(M,\mathbf{E})}
\leq C_s\|\mathbf{u}\|_{H^s(M,\mathbf{E})}.
\label{pseudo:eq-estimacion-resolvente-sectorial}
\end{equation}
The constant $C_s$ depends on $s$, the sector, finitely many parameter seminorms of the symbol, and the choices of a finite atlas, but not on $\mathbf{u}$ or $\lambda$ in the indicated region; the estimate holds for every $\mathbf{u}\in H^s(M,\mathbf{E})$.
\end{proposition}

\begin{proof}
The symbolic composition formula gives
\[
(p\#q)(x,\xi,\lambda)
\sim
\sum_{\alpha\in\mathbb N_0^n}
\frac{i^{-|\alpha|}}{\alpha!}
(D_\xi^\alpha p)(D_x^\alpha q).
\]
Transfer of Euclidean weights in the composition estimates is governed by Peetre's inequality in Lemma~\ref{lem:desigualdad-peetre-peso-japones}; the weight $\Omega_m(\xi,\lambda)$ also preserves the resolvent factor recorded separately below. The term of degree zero in the total weight is $p_mq_{-m}$; imposing $p_mq_{-m}=I$ gives $q_{-m}=p_m^{-1}$. The parameter-ellipticity estimate proves \eqref{pseudo:eq-cotas-componentes-resolvente} for $j=0$, $\alpha=\beta=0$, and $k=0$. If $\mathscr D$ is one of the derivatives $D_{x^r}$, $D_{\xi_r}$, or $\partial_\lambda$, differentiating $p_mq_{-m}=I$ gives the exact identity
\begin{equation}
 \mathscr Dq_{-m}=-q_{-m}(\mathscr Dp_m)q_{-m}.
 \label{pseudo:eq-derivada-inversa-principal-parametro}
\end{equation}
For $\mathscr D=\partial_{\bar\lambda}$, the right-hand side is zero. To order the derivatives together, write $\mathscr D^\nu=D_x^\alpha D_\xi^\beta\partial_\lambda^h$ and use the componentwise order on triple multi-indices. Applying the Leibniz rule to $\mathscr D^\nu(p_mq_{-m})=0$ and then multiplying on the left by $q_{-m}$ gives the exact recurrence
\begin{equation}
 \mathscr D^\nu q_{-m}
 =-q_{-m}\sum_{0<\delta\leq\nu}
 \binom{\nu}{\delta}
 (\mathscr D^\delta p_m)
 (\mathscr D^{\nu-\delta}q_{-m}).
 \label{pseudo:eq-recurrencia-derivadas-inversa-principal-parametro}
\end{equation}
This formula closes an induction on $|\nu|$: in each summand, the derivative of $q_{-m}$ in the last factor has strictly smaller order. If the $\lambda$ component of $\delta$ is zero, the symbol bounds on $p_m$ and the induction hypothesis bound this summand by
\[
 C\Omega_m^{-1}
 \langle\xi\rangle^{m-|\delta_\xi|}
 \langle\xi\rangle^{-|\beta|+|\delta_\xi|}
 \Omega_m^{-1-h}
 \leq
 C\langle\xi\rangle^{-|\beta|}\Omega_m^{-1-h}.
\]
If that component is one, the summand can be nonzero only when $\mathscr D^\delta p_m=I$; the induction hypothesis then gives $C\langle\xi\rangle^{-|\beta|}\Omega_m^{-1-h}$ directly. For two or more derivatives in $\lambda$, the differentiated factor of $p_m$ is zero. The sum in \eqref{pseudo:eq-recurrencia-derivadas-inversa-principal-parametro} is finite, so
\[
 \|D_x^\alpha D_\xi^\beta\partial_\lambda^kq_{-m}\|
 \leq C_{K,\alpha,\beta,k}
 \langle\xi\rangle^{-|\beta|}\Omega_m^{-1-k}.
\]
This establishes the initial case, including holomorphy.

Now suppose the components $q_{-m},\ldots,q_{-m-j+1}$ have been constructed and their homogeneity, holomorphy, and \eqref{pseudo:eq-cotas-componentes-resolvente} proved. The coefficient of degree $-j$ in $p\# q$ is
\[
 p_mq_{-m-j}
 +\sum_{\substack{k+\ell+|\gamma|=j\\0\leq\ell<j}}
 \frac{i^{-|\gamma|}}{\gamma!}
 (D_\xi^\gamma p_{m-k})(D_x^\gamma q_{-m-\ell}).
\]
Setting it equal to zero and multiplying on the left by $q_{-m}$ gives exactly \eqref{pseudo:eq-recurrencia-parametrix-resolvente}; thus the new component is unique. Under dilation $(\xi,\lambda)\mapsto(t\xi,t^m\lambda)$, the product inside the sum has degree
\[
 (m-k-|\gamma|)+(-m-\ell)=-(k+\ell+|\gamma|)=-j.
\]
The outer factor $q_{-m}$ lowers the degree to $-m-j$. Every factor in the recurrence is holomorphic in $\lambda$, so the new component is also holomorphic.

To prove its estimates, now write the symbol index as $\kappa$ and the number of derivatives in $\lambda$ as $h$. Distributing $D_x^\alpha D_\xi^\beta\partial_\lambda^h$ among the three factors of a summand with $\kappa+\ell+|\gamma|=j$ gives Leibniz terms of the form, apart from their multinomial coefficient,
\[
 \bigl(D_x^{\alpha_0}D_\xi^{\beta_0}
       \partial_\lambda^{h_0}q_{-m}\bigr)
 \bigl(D_x^{\alpha_1}D_\xi^{\beta_1+\gamma}
       \partial_\lambda^{h_1}p_{m-\kappa}\bigr)
 \bigl(D_x^{\alpha_2+\gamma}D_\xi^{\beta_2}
       \partial_\lambda^{h_2}q_{-m-\ell}\bigr),
\]
where $\alpha_0+\alpha_1+\alpha_2=\alpha$, $\beta_0+\beta_1+\beta_2=\beta$, and $h_0+h_1+h_2=h$. If $h_1=0$, the induction hypotheses bound its norm by
\[
 C\langle\xi\rangle^{-j-|\beta|+m}
 \Omega_m^{-2-h}
 \leq C\langle\xi\rangle^{-j-|\beta|}\Omega_m^{-1-h}.
\]
If $h_1=1$, a nonzero middle factor would require $\kappa=0$,
$\gamma=0$, and $\alpha_1=\beta_1=0$. The recurrence would then give
$j=\ell$, contrary to $\ell<j$. Thus every term with $h_1=1$ vanishes. If $h_1\geq2$, the term is zero. There are finitely many Leibniz distributions and finitely many triples $(\kappa,\ell,\gamma)$ for each $j$, so summing proves \eqref{pseudo:eq-cotas-componentes-resolvente} for the new component. This completes the induction.

For two-sided compatibility, the equation $q\#p=I$ gives, at degree $-j$, the equation \[
 q_{-m-j}p_m
 +\sum_{\substack{k+\ell+|\alpha|=j\\0\leq\ell<j}}
 \frac{i^{-|\alpha|}}{\alpha!}
 (D_\xi^\alpha q_{-m-\ell})(D_x^\alpha p_{m-k})=0,
\], and multiplying on the right by $q_{-m}$ gives \eqref{pseudo:eq-recurrencia-derecha-parametrix-resolvente}. To check that this does not define different components, let $q^{\mathrm d}$ be the formal series determined by $p\#q^{\mathrm d}=I$ and $q^{\mathrm i}$ the one determined by $q^{\mathrm i}\#p=I$. Formal associativity, whose coefficient at each degree contains only finite sums, gives
\[
 q^{\mathrm i}=q^{\mathrm i}\#(p\#q^{\mathrm d})
 =(q^{\mathrm i}\#p)\#q^{\mathrm d}=q^{\mathrm d}.
\]
Thus the two recurrences determine the same series term by term. This entire construction is at high frequency: $a_{m-k}$ and $q_{-m-j}$ are evaluated only for $\xi\neq0$.

To complete the construction on the axis $\xi=0$, use the full symbol \[
p(x,\xi,\lambda):=\lambda I-a(x,\xi),
\], which is smooth for every $\xi$. There exists $R>0$ such that $p$ is pointwise invertible when $\widehat\rho_m\geq R$. To see this quantitatively, fix $C_0>0$. In the cone $|\lambda|\leq C_0\langle\xi\rangle^m$, the principal estimate and the bound $\|a-a_m\|\leq C_1\langle\xi\rangle^{m-1}$ imply
\[
 \|p_m^{-1}(a-a_m)\|\leq C_2\langle\xi\rangle^{-1}.
\]
If $\langle\xi\rangle\geq2C_2$, the matrix geometric series inverts $p=p_m-(a-a_m)=p_m[I-p_m^{-1}(a-a_m)]$. In the complementary cone, choose $C_0\geq2C_3$, where $\|a(x,\xi)\|\leq C_3\langle\xi\rangle^m$; then \[
 \|\lambda^{-1}a(x,\xi)\|\leq\frac{1}{2}
 \quad\text{if}\quad |\lambda|\geq C_0\langle\xi\rangle^m,
\] and $p=\lambda[I-\lambda^{-1}a]$ is invertible. Comparison \eqref{pseudo:eq-comparacion-pesos-parametro} allows a single $R$ to be chosen for both cones.

Let $\chi\in C^\infty(\mathbb R)$ vanish on $(-\infty,1]$ and equal one on $[2,\infty)$. Extending by zero where the inverse is not used, define
\[
b_0(x,\xi,\lambda)
:=
\chi\!\left(\frac{\widehat\rho_m(\xi,\lambda)}{R}\right)
p(x,\xi,\lambda)^{-1}.
\]
The inverse-derivative formula \eqref{pseudo:eq-derivada-inversa-principal-parametro}, now applied to $p^{-1}$, the symbol bounds on $a$, and \eqref{pseudo:eq-derivadas-corte-rho-parametro} show, derivative by derivative, that $b_0$ is a global family of resolvent type $0$. It is not holomorphic, since the cutoff depends on $|\lambda|^2$. For \[
e:=I-p\# b_0
\], we have
$e$
of defect type $1$: the contribution $I-\chi$ has jointly compact support in $(\xi,\lambda)$, and every remaining term in the symbolic calculus contains at least one derivative in $\xi$. Outside that compact support, we therefore have the sharper bound
\begin{equation}
\left\|
D_x^\alpha D_\xi^\beta
\partial_\lambda^k\partial_{\bar\lambda}^\ell e
\right\|
\leq
C_{\alpha,\beta,k,\ell}
\langle\xi\rangle^{m-1-|\beta|}
\Omega_m(\xi,\lambda)^{-1-k-\ell}.
\label{pseudo:eq-defecto-parametrico-reforzado}
\end{equation}
That is, besides losing one power of $\langle\xi\rangle$, the defect contains the factor $\langle\xi\rangle^{m-1}\Omega_m^{-1}$. To obtain powers of the defect, start from $e^{\#(N+1)}=e^{\#N}\#e$. If the following bound holds for $N$, the composition formula distributes derivatives between $e^{\#N}$ and $e$; the product of the weight factors is
\[
 \langle\xi\rangle^{N(m-1)}\Omega_m^{-N}
 \langle\xi\rangle^{m-1}\Omega_m^{-1}
 =\langle\xi\rangle^{(N+1)(m-1)}\Omega_m^{-(N+1)}.
\]
Each derivative in $\xi$ lowers the power of $\langle\xi\rangle$, and each Wirtinger derivative lowers that of $\Omega_m$; the Leibniz rule produces a finite sum. The case $N=1$ is \eqref{pseudo:eq-defecto-parametrico-reforzado}, so induction gives, modulo a family of jointly compact support,
\[
\left\|D_x^\alpha D_\xi^\beta
\partial_\lambda^k\partial_{\bar\lambda}^\ell e^{\# N}\right\|
\leq
C_{N,\alpha,\beta,k,\ell}
\langle\xi\rangle^{N(m-1)-|\beta|}
\Omega_m^{-N-k-\ell}.
\]
Given $P,Q\in\mathbb N_0$, choose an integer $N\geq\displaystyle\max\{Q,\,P+mQ+1\}$. Since $\Omega_m\geq1+|\lambda|$, $\Omega_m\geq\langle\xi\rangle^m$, and $N-Q\geq0$, we have
\[
\Omega_m^{-N-k-\ell}
\leq(1+|\lambda|)^{-Q-k-\ell}
\langle\xi\rangle^{-m(N-Q)}.
\]
Thus the preceding inequality implies
\begin{equation}
 \left\|D_x^\alpha D_\xi^\beta
 \partial_\lambda^k\partial_{\bar\lambda}^\ell e^{\# N}\right\|
 \leq C_{P,Q,N,\alpha,\beta,k,\ell}
 \langle\xi\rangle^{-P-|\beta|}(1+|\lambda|)^{-Q-k-\ell}.
 \label{pseudo:eq-potencias-defecto-regularizantes-parametro}
\end{equation}
Indeed, after extracting $(1+|\lambda|)^{-Q-k-\ell}$, the remaining power is $\langle\xi\rangle^{N(m-1)-m(N-Q)}=
\langle\xi\rangle^{-N+mQ}$, which does not exceed $\langle\xi\rangle^{-P}$. This is quantified simultaneous smoothing. In particular, $c_k:=b_0\#e^{\# k}$ has resolvent type $k$, and for every $P,Q$ and each finite family of derivatives there exists an integer $J=J(P,Q,\alpha,\beta,k,\ell)$ such that terms of index $j\geq J$ satisfy the filtered simultaneous-smoothing hypothesis. Apply the second part of Lemma~\ref{pseudo:lem-suma-asintotica-parametro-debil}, with the joint cutoffs $\chi(\varepsilon_k\widehat\rho_m)$, to sum the formal series
\[
\widetilde q
\sim
b_0\#\sum_{k=0}^{\infty}e^{\# k}
\]
and obtain the nonholomorphic auxiliary family $\widetilde q$. If \[
 T_N:=b_0\#\sum_{k=0}^{N-1}e^{\# k},
\], the finite geometric identity gives
\begin{equation}
 p\#T_N=(I-e)\#\sum_{k=0}^{N-1}e^{\# k}=I-e^{\#N}.
 \label{pseudo:eq-identidad-neumann-finita-parametro}
\end{equation}
In the algebra of full symbols, possible corrections from the composition formula already belong to the ideal of families smoothing with parameter. Thus
\[
 s_1:=I-p\#\widetilde q
 =e^{\#N}-p\#(\widetilde q-T_N).
\]
Given $P,Q$ and a finite family of derivatives, choose $N$ as in \eqref{pseudo:eq-potencias-defecto-regularizantes-parametro}. The first part of the last expression satisfies the required bound there, and the diagonal choice of cutoffs in the lemma gives the same bound for the second. Thus
\[
 \|D_x^\alpha D_\xi^\beta\partial_\lambda^k
 \partial_{\bar\lambda}^{\ell}s_1\|
 \leq C_{P,Q,\alpha,\beta,k,\ell}
 \langle\xi\rangle^{-P}(1+|\lambda|)^{-Q}.
\]
The left construction provides $\widetilde q_{\mathrm i}\#p=I-s_2$ with the same estimate. The identity \[
 \widetilde q-\widetilde q_{\mathrm i}
 =s_2\#\widetilde q-\widetilde q_{\mathrm i}\#s_1
\] shows that both families differ by a family smoothing with parameter. Replacing one by the other preserves its expansion and yields a two-sided inverse modulo that ideal. Expanding $p$ and $\widetilde q$ for large $\|\xi\|$, uniqueness at each degree recovers exactly the components $q_{-m-j}$ of \eqref{pseudo:eq-recurrencia-parametrix-resolvente}.

The cutoffs in this full sum depend on $\widehat\rho_m$ and therefore need not be holomorphic. We will not use $\widetilde q$ as though it were: contour deformations will be performed only on each recursive component; the exact resolvent, although holomorphic on its resolvent set, will always be integrated over the fixed global contour.

Quantizing $\widetilde q$ in charts, inserting partitions of unity, and summing gives $Q(\lambda)$ with
\[
(\lambda I-A)Q(\lambda)=I-S_1(\lambda),
\qquad
Q(\lambda)(\lambda I-A)=I-S_2(\lambda),
\]
where $S_1$ and $S_2$ are smoothing. Relative to the measure $d\lambda_{\mathbf{g}}(y)$, their kernels satisfy
\[
\mathbf{K}_{S_j}(\lambda;\cdot,\cdot)
\in C^\infty\!\left(
M\times M;
\operatorname{Hom}(\operatorname{pr}_2^*\mathbf{E},\operatorname{pr}_1^*\mathbf{E})
\right).
\]
Equip this bundle with the connection induced by the Hermitian connection of $\mathbf{E}$ and the tensor norms induced by $\mathbf{g}$ and $\mathbf{h}_{\mathbf{E}}$. For all $Q,k,\ell,a,b\in\mathbb N_0$ there exists $C_{Q,k,\ell,a,b}>0$ such that
\begin{equation}
 \sup_{x,y\in M}
 \left|
 \nabla_x^a\nabla_y^b\partial_\lambda^k
 \partial_{\bar\lambda}^{\ell}\mathbf{K}_{S_j}(\lambda;x,y)
 \right|_{\mathbf{g},\mathbf{h}_{\mathbf{E}}}
 \leq C_{Q,k,\ell,a,b}(1+|\lambda|)^{-Q}.
 \label{pseudo:eq-residuos-resolvente-decaimiento-cuantificado}
\end{equation}
Let $V=\operatorname{vol}_{\mathbf{g}}(M)$. Proposition~\ref{prop:criterio-schur-operadores-integrales}, applied fiberwise with $X=Y=M$ and $\mu=\nu=\lambda_{\mathbf{g}}$, together with \eqref{pseudo:eq-residuos-resolvente-decaimiento-cuantificado} for $\alpha=\beta=k=\ell=0$, gives, for each $Q\in\mathbb N$,
\begin{equation}
 \|S_j(\lambda)\|_{\mathcal L(L^2(M,\mathbf{E}))}
 \leq V\sup_{x,y\in M}
 \|\mathbf{K}_{S_j}(\lambda;x,y)\|_{
 \mathcal L((\mathbf{E}_y,\mathbf{h}_{\mathbf{E}}(y)),(\mathbf{E}_x,\mathbf{h}_{\mathbf{E}}(x)))}
 \leq VC_Q(1+|\lambda|)^{-Q}.
 \label{pseudo:eq-pequenez-residuos-neumann}
\end{equation}
Choose $R_1$ so that the last expression is at most $\displaystyle\frac{1}{2}$ for $j=1,2$, $\lambda\in
\Lambda^{\mathrm p}_{\theta,\varepsilon_{\mathrm p}}$, and $|\lambda|\geq R_1$. The series
$\displaystyle \sum_{\nu=0}^{\infty}S_j(\lambda)^\nu$
then converge in $\mathcal L(L^2(M,\mathbf{E}))$, their sums have norm at most two, and the two parametrix identities provide a right and a left inverse of $\lambda I-A$, respectively. They coincide and define the global resolvent. Thus
\begin{equation}
 \{\lambda\in\Lambda^{\mathrm p}_{\theta,\varepsilon_{\mathrm p}}\mid
 |\lambda|\geq R_1\}\subseteq\rho(A).
 \label{pseudo:eq-sector-principal-resolvente-grande}
\end{equation}

Narrowing the sector is now justified. The set \[
 \mathbf{F}:=\sigma(A)\cap\{\lambda\mid |\lambda|\leq R_1\}
\] is finite: admissibility of the cut implies $L_\theta\subseteq\rho(A)$, so Proposition~\ref{pseudo:prop-resolvente-compacto-espectro-discreto} applies. Moreover, $0\notin \mathbf{F}$ and $\mathbf{F}\cap L_\theta=\varnothing$ by the cut hypotheses. Finiteness of $\mathbf{F}$ gives $0<\varepsilon\leq\varepsilon_{\mathrm p}$ such that $\mathbf{F}\cap\Lambda_{\theta,\varepsilon}=\varnothing$. Together with \eqref{pseudo:eq-sector-principal-resolvente-grande}, this proves $\Lambda_{\theta,\varepsilon}\subseteq\rho(A)$. Fix $R_0\geq\displaystyle\max\{R_1,1\}$.

For $|\lambda|\geq R_0$, the algebraic identity
\[
(\lambda I-A)^{-1}
=Q(\lambda)(I-S_1(\lambda))^{-1}
=(I-S_2(\lambda))^{-1}Q(\lambda)
\]
shows that the result is the exact resolvent; its holomorphy follows from holomorphy of $\lambda\mapsto\lambda I-A$ and uniqueness of the inverse, not from the cutoffs used in summation. Explicitly, the resolvent identity gives
\[
 \frac{(\lambda+h-A)^{-1}-(\lambda-A)^{-1}}{h}
 =-(\lambda+h-A)^{-1}(\lambda-A)^{-1},
\]
and the right-hand side converges in operator norm to $-(\lambda-A)^{-2}$ as $h\to0$. This proves operator holomorphy of the exact resolvent. The identity \[
 (\lambda I-A)^{-1}-Q(\lambda)
 =Q(\lambda)\bigl[(I-S_1(\lambda))^{-1}-I\bigr]
\] and its left version, together with \eqref{pseudo:eq-residuos-resolvente-decaimiento-cuantificado}, show that the difference satisfies the same estimate for each $Q$ and each finite family of kernel derivatives.

On the part $|\lambda|\leq R_0$ of the sector, the resolvent exists by the choice just justified. Since $0\in\rho(A)$, the closure of this part after adding the vertex is a compact subset of $\rho(A)$. The two parametrix identities give
\[
(\lambda I-A)^{-1}-Q(\lambda)
=(\lambda I-A)^{-1}S_1(\lambda)
=S_2(\lambda)(\lambda I-A)^{-1}.
\]
Elliptic regularity applied on both sides shows that this difference is smoothing; on each compact subset of the sector, all its kernel seminorms are uniform. The difference may fail to be holomorphic because $Q$ contains cutoffs depending on $|\lambda|$.

In a fixed chart and quantization, extracting the full symbol from the kernel is a continuous linear operation. Applied directly to the exact resolvent, it produces a symbol $r(x,\xi,\lambda)$ holomorphic throughout the sector. For large parameters, $r-\widetilde q$ is smoothing with parameter, and on the compact part it is uniformly smoothing; thus $r$ retains the high-frequency expansion formed by the $q_{-m-j}$. Bounds on its derivatives in $\lambda$ follow from Cauchy inequalities in a slightly smaller sector. Finally, the weak estimates of resolvent type $0$ give the first norm in \eqref{pseudo:eq-estimacion-resolvente-sectorial}; using $(\lambda I-A)(\lambda I-A)^{-1}=I$ gives the second. The construction is matrix-valued from the first step, so it is unchanged for bundle sections.
\end{proof}

\subsection{Construction of complex powers}

The parameter-dependent parametrix gives the high-frequency behavior of the resolvent, while the exact resolvent retains holomorphy and spectral information. Integrating them against $\lambda^z$ over a contour around the cut combines both properties in an operator of complex order.

Let $\Gamma_\theta$ be a keyhole contour contained in $\rho(A)$, with two infinite branches in the spectrum-free sector and a small arc around zero. It is the limit, as the outer arc is removed, of positively oriented closed contours surrounding the spectrum of $A$ in the cut plane. With this orientation and the resolvent $(\lambda I-A)^{-1}$, the scalar residue at an eigenvalue $\mu$ is $+\mu^z$. Invertibility prevents the small arc from meeting the spectrum, and estimate \[
 \|(\lambda I-A)^{-1}\|_{\mathcal L(L^2(M,\mathbf{E}))}
 \leq C_0|\lambda|^{-1},
 \qquad
 \lambda\in\Lambda_{\theta,\varepsilon},\quad |\lambda|\geq R_0,
\] makes the following integral converge in operator norm when $\operatorname{Re}z<0$.

\begin{definition}[Complex powers]
\label{pseudo:def-potencias-complejas}
\index{complex power!of an elliptic operator}
Let $(M,\mathbf{g})$ be a smooth closed Riemannian manifold; let $\mathbf{E}\longrightarrow M$ be a complex vector bundle of finite rank equipped with a Hermitian bundle metric; let $A\in\Psi^m_{\mathrm{cl}}(M;\mathbf{E},\mathbf{E})$ be elliptic and invertible, with $m>0$; and let $\theta$ be an admissible spectral cut. Let $\Gamma_\theta\subset\rho(A)$ be the positively oriented keyhole contour described above, and use the branch $\log_\theta$ in $\lambda^z$. For $\operatorname{Re}z<0$, define
\begin{equation}
A^z
:=
\frac{1}{2\pi i}
\int_{\Gamma_\theta}
\lambda^z(\lambda I-A)^{-1}\,d\lambda.
\label{pseudo:eq-def-potencia-negativa}
\end{equation}
For arbitrary $z\in\mathbb C$, choose $k\in\mathbb N$ with $\operatorname{Re}(z-k)<0$ and set
\begin{equation}
A^z:=A^kA^{z-k}
\quad\text{on }\Gamma(M,\mathbf{E}).
\label{pseudo:eq-extension-potencias}
\end{equation}
\end{definition}

The resolvent identity
\[
(\lambda I-A)^{-1}-(\mu I-A)^{-1}
=(\mu-\lambda)(\lambda I-A)^{-1}(\mu I-A)^{-1}
\]
and two successive applications of Cauchy's formula first prove, for negative real parts, that
\[
A^zA^w=A^{z+w}.
\]

Let us check that the integral with $z=-1$ is the usual inverse of $A$. Indeed, the contour calculus for rational functions follows directly from the resolvent identity and Cauchy's formula; since $0\in\rho(A)$, applying it to $r(\lambda)=\lambda^{-1}$ gives
\begin{equation}
 \frac{1}{2\pi i}
 \int_{\Gamma_\theta}
 \lambda^{-1}(\lambda I-A)^{-1}\,d\lambda
 =A^{-1}.
 \label{pseudo:eq-contorno-inverso-usual}
\end{equation}
Let us give the analytic justification of this identity. Truncate the keyhole contour and apply Cauchy's formula before removing the truncation. On $\Gamma(M,\mathbf{E})$, use
\[
 A(\lambda I-A)^{-1}
 =\lambda(\lambda I-A)^{-1}-I,
 \qquad
 (\lambda I-A)^{-1}A
 =\lambda(\lambda I-A)^{-1}-I.
\]
Moreover, for $\mathbf{u}\in\Gamma(M,\mathbf{E})$,
\[
 (\lambda I-A)^{-1}\mathbf{u}-\lambda^{-1}\mathbf{u}
 =\lambda^{-1}(\lambda I-A)^{-1}A\mathbf{u}.
\]
The resolvent estimate makes the remaining boundary terms tend to zero as the outer arc is removed, while the inner arc remains inside $\rho(A)$ because $0\in\rho(A)$. The scalar Cauchy formula then gives, respectively,
\[
 AA^{-1}\mathbf{u}=\mathbf{u},
 \qquad
 A^{-1}A\mathbf{u}=\mathbf{u},
 \qquad \mathbf{u}\in\Gamma(M,\mathbf{E}),
\]
proving \eqref{pseudo:eq-contorno-inverso-usual} without identifying the integral with the inverse in advance.

We can now check that \eqref{pseudo:eq-extension-potencias} is independent of $k$. If $\operatorname{Re}\zeta<-1$, the integrals converge after multiplication by $A$, and the preceding identities give
\[
 AA^\zeta=A^\zeta A=A^{\zeta+1};
\]
the additional term $\displaystyle \frac{1}{2\pi i}\int_{\Gamma_\theta}\lambda^\zeta\,d\lambda$ is zero by Cauchy's formula and decay at infinity. If $k$ is admissible and $\zeta:=z-k-1$, then $\operatorname{Re}\zeta<-1$ and
\[
 A^{k+1}A^{z-k-1}
 =A^k(AA^\zeta)
 =A^kA^{z-k}.
\]
By induction, any two admissible choices produce the same operator. The identity \eqref{pseudo:eq-contorno-inverso-usual} gives $A^0=I$ and shows that integer powers are the usual ones.

Finally, for $z,w\in\mathbb C$ choose $k,\ell\in\mathbb N$ with $\operatorname{Re}(z-k)<0$ and $\operatorname{Re}(w-\ell)<0$. Resolvents commute with $A$ on $\Gamma(M,\mathbf{E})$, so negative powers commute with integer powers. Using the law already proved in the negative half-plane,
\[
 A^zA^w
 =A^{k+\ell}A^{z-k}A^{w-\ell}
 =A^{k+\ell}A^{z+w-k-\ell}
 =A^{z+w}.
\]
Invertibility is precisely what permits this group law for all $z\in\mathbb C$; without it, negative powers require a separate convention on the kernel.

For $\mu\in\mathbb C$, the notation $\Psi^\mu_{\mathrm{cl}}(M;\mathbf{E},\mathbf{E})$ means that the symbol has homogeneous components of degrees $\mu-j$ and satisfies estimates of order $\operatorname{Re}\mu$. For $t>0$, interpret $t^\mu=e^{\mu\log t}$ using the real logarithm.

\begin{definition}[Holomorphy for a family of variable order]
\label{pseudo:def-holomorfia-familia-orden-variable}
Let $(M,\mathbf{g})$ be a smooth closed Riemannian manifold, let $\mathbf{E}\longrightarrow M$ be a complex vector bundle of finite rank equipped with a Hermitian bundle metric, let $m>0$, and let $G\subseteq\mathbb C$ be open. A family $B(z)\in\Psi^{mz}_{\mathrm{cl}}(M;\mathbf{E},\mathbf{E})$ is \textbf{holomorphic of variable order $mz$} if the following three levels hold.

\emph{Homogeneous components.} In every chart and for $\xi\neq0$, the components $b_{mz-j}(x,\xi;z)$ are holomorphic in $G$ (entire when $G=\mathbb C$), and for each $K\Subset G$, compact subset $K_x$ of the chart, $h\in\mathbb N_0$, and $\alpha,\beta$, there exists $C>0$ such that
\begin{equation}
 \|D_x^\alpha D_\xi^\beta\partial_z^h
 b_{mz-j}(x,\xi;z)\|
 \leq C\|\xi\|^{m\operatorname{Re}z-j-|\beta|}
 (1+|\log\|\xi\||)^h,
 \quad \|\xi\|\geq1,\quad z\in K.
 \label{pseudo:eq-holomorfia-componentes-potencias}
\end{equation}

\emph{Full symbol.} With a fixed homogeneous cutoff $\chi(\xi)$, the full symbol $b(x,\xi;z)$ is holomorphic in $z$ and, for each $N,h\in\mathbb N_0$, $\delta>0$, $K\Subset G$, and $\alpha,\beta$, there exists $C>0$ such that
\begin{equation}
 \left\|D_x^\alpha D_\xi^\beta\partial_z^h
 \left(b-\sum_{j=0}^{N-1}\chi b_{mz-j}\right)(x,\xi;z)\right\|
 \leq C\langle\xi\rangle^{m\operatorname{Re}z-N-|\beta|+\delta},
 \qquad z\in K.
 \label{pseudo:eq-holomorfia-simbolo-completo-potencias}
\end{equation}
The arbitrary loss $\delta$ absorbs the powers of $\log\langle\xi\rangle$ arising from differentiation in $z$. Relative to $d\lambda_{\mathbf{g}}(y)$, the smoothing remainder has a kernel \[
\mathbf{K}_R(z;\cdot,\cdot)
\in C^\infty\!\left(
M\times M;
\operatorname{Hom}(\operatorname{pr}_2^*\mathbf{E},\operatorname{pr}_1^*\mathbf{E})
\right)
\] holomorphic in $z$. With the induced connection and tensor norm, for all $h,a,b\in\mathbb N_0$,
\[
 \sup_{\substack{z\in K\\x,y\in M}}
 \left|\partial_z^h\nabla_x^a\nabla_y^b\mathbf{K}_R(z;x,y)\right|_{\mathbf{g},\mathbf{h}_{\mathbf{E}}}
 \leq C_{K,h,a,b}.
\]

\emph{Operator level.} The map $z\mapsto B(z)$ is holomorphic from $G$ into $\mathcal L_b(\Gamma(M,\mathbf{E}),\Gamma(M,\mathbf{E}))$, where $\Gamma(M,\mathbf{E})$ carries its Fréchet topology and the subscript $b$ denotes uniform convergence on bounded subsets. In the presence of the two preceding symbol levels, this is equivalent to the following fixed-Sobolev formulation. If $\displaystyle \tau_K:=m\displaystyle\max_{z\in K}\operatorname{Re}z$, then, for every $r\in\mathbb R$, $h\in\mathbb N_0$, and $\delta>0$,
\begin{equation}
 \sup_{z\in K}
 \|\partial_z^hB(z)\|_{\mathcal L(
 H^r(M,\mathbf{E}),H^{r-\tau_K-\delta}(M,\mathbf{E}))}
 \leq C_{K,r,h,\delta},
 \label{pseudo:eq-holomorfia-operatorial-potencias}
\end{equation}
and $z\mapsto B(z)$ is holomorphic on an open neighborhood of $K$ contained in $G$, with values in that bounded-operator space. Choose this neighborhood so that $m\operatorname{Re}z<\tau_K+\displaystyle\frac{\delta}{2}$; the remaining margin absorbs the logarithmic powers of the derivatives. The entire family over $G$ need not take values in a single Sobolev pair. The fixed order $\tau_K+\delta$ allows operators whose real orders vary with $z$ to be compared in the same Banach space.
\end{definition}

\begin{theorem}[Complex powers of an elliptic operator]
\label{pseudo:teo-potencias-complejas-clasicas}
Let $(M,\mathbf{g})$ be a smooth closed Riemannian manifold; let $\mathbf{E}\longrightarrow M$ be a complex vector bundle of finite rank equipped with a Hermitian bundle metric $\mathbf{h}_{\mathbf{E}}$ and a compatible connection $\nabla^{\mathbf{E}}$; let $A\in\Psi^m_{\mathrm{cl}}(M;\mathbf{E},\mathbf{E})$ be elliptic and invertible, with $m>0$ and principal Fourier symbol $a_m$; and fix an admissible spectral cut $\theta$. For every $z\in\mathbb C$,
\[
A^z\in\Psi^{mz}_{\mathrm{cl}}(M;\mathbf{E},\mathbf{E}).
\]
The family is holomorphic in the sense of Definition~\ref{pseudo:def-holomorfia-familia-orden-variable} and satisfies
\begin{equation}
A^zA^w=A^{z+w},
\qquad
\boldsymbol{\sigma}_{mz}^\Psi(A^z)(x,\xi)=a_m(x,\xi)^z,
\label{pseudo:eq-ley-simbolo-potencias}
\end{equation}
where the matrix power on the right uses the same spectral cut. Moreover,
\[
A^z\colon H^r(M,\mathbf{E})\longrightarrow H^{r-m\operatorname{Re}z}(M,\mathbf{E})
\]
is continuous for every $r\in\mathbb R$ when initially considered on smooth sections.
\end{theorem}

\begin{proof}
Begin with $\operatorname{Re}z<0$. By Proposition~\ref{pseudo:prop-parametrix-resolvente}, the resolvent symbol is \[
r(x,\xi,\lambda)
\underset{\|\xi\|\to\infty}{\sim}
\sum_{j=0}^\infty q_{-m-j}(x,\xi,\lambda).
\]. Extraction of the full symbol is linear and continuous, so it commutes with the convergent integral defining $A^z$. Thus a symbol of $A^z$ is
\[
b(x,\xi;z)
=
\frac{1}{2\pi i}
\int_{\Gamma_\theta}\lambda^z r(x,\xi,\lambda)\,d\lambda.
\]
The following deformation is performed exclusively for homogeneous components, not for $r$. Indeed, compactness of the cosphere bundle, invertibility of $a_m$, and the principal angle give a normalized contour uniformly avoiding $\operatorname{spec}a_m(x,\omega)$ for every $\|\omega\|=1$. Let $\Gamma_{\theta,\xi}$ denote its image under $\lambda=\|\xi\|^m\mu$. Its branches lie in the principal sector, and its curved part has radius $c\|\xi\|^m$, with $c>0$ smaller than the modulus of every eigenvalue of $a_m(x,\omega)$. For large $\|\xi\|$, the fixed contour $\Gamma_\theta$ and $\Gamma_{\theta,\xi}$ are homotopic in the set where $\lambda I-a_m(x,\xi)$ is invertible. This homotopy may cross global eigenvalues of $A$ and therefore does not apply to the exact symbol $r$; it does apply to each $q_{-m-j}$, whose only singularities in $\lambda$ arise from the principal resolvent. For $\xi\neq0$, define
\begin{equation}
b_{mz-j}(x,\xi;z)
:=
\frac{1}{2\pi i}
\int_{\Gamma_{\theta,\xi}}
\lambda^zq_{-m-j}(x,\xi,\lambda)\,d\lambda.
\label{pseudo:eq-componentes-simbolo-potencia}
\end{equation}
Parameter estimates allow differentiation under the integral. Indeed, if $K\Subset\{z\mid\operatorname{Re}z<0\}$ and $\displaystyle \eta:=-\displaystyle\max_{z\in K}\operatorname{Re}z>0$, then on the branches where $|\lambda|\geq1$,
\begin{equation}
 \bigl\|\partial_z^h\bigl(\lambda^z
 D_x^\alpha D_\xi^\beta q_{-m-j}\bigr)\bigr\|
 \leq C_{K,h,\alpha,\beta,j}
 (1+\log|\lambda|)^h|\lambda|^{-1-\eta}
 \langle\xi\rangle^{-j-|\beta|}.
 \label{pseudo:eq-dominante-integral-componentes-potencia}
\end{equation}
The function $(1+\log r)^hr^{-1-\eta}$ is integrable in $[1,\infty)$; after normalization by $\lambda=\|\xi\|^m\mu$, the arc belongs to a compact set uniformly avoiding the principal spectrum. This justifies, in the symbol topology, both the integral and all its derivatives in $z,x$ and $\xi$. Replacing $\xi$ by $t\xi$ and making the change $\lambda=t^m\mu$, homogeneity of the symbolic resolvent gives
\[
\begin{aligned}
b_{mz-j}(x,t\xi;z)
&=
\frac{1}{2\pi i}
\int_{\Gamma_{\theta,\xi}}
(t^m\mu)^z
t^{-m-j}q_{-m-j}(x,\xi,\mu)t^m\,d\mu\\
&=t^{mz-j}b_{mz-j}(x,\xi;z).
\end{aligned}
\]
Thus the terms in \eqref{pseudo:eq-componentes-simbolo-potencia} have exactly the required degrees.

It remains to control the remainder without moving the global resolvent contour. Fix a cutoff $\chi_0(\xi)$ vanishing near zero and equal to one beyond a frequency for which the preceding principal homotopy is valid. For $L\in\mathbb N$, set, on the fixed contour,
\[
 R_L(x,\xi,\lambda)
 :=
 r(x,\xi,\lambda)
 -\sum_{j=0}^{L-1}\chi_0(\xi)q_{-m-j}(x,\xi,\lambda).
\]
The parameter-dependent expansion of Proposition~\ref{pseudo:prop-parametrix-resolvente} and uniform smoothing of $r-\widetilde q$ give, on the two infinite branches,
\begin{equation}
 \|D_x^\alpha D_\xi^\beta R_L(x,\xi,\lambda)\|
 \leq C_{\alpha,\beta,L}
 \langle\xi\rangle^{-L-|\beta|}
 \Omega_m(\xi,\lambda)^{-1}.
 \label{pseudo:eq-resto-resolvente-contorno-fijo}
\end{equation}
On the fixed curved part of the contour, the same inequality for large $\|\xi\|$ is the ordinary estimate of order $-m-L$; the region of bounded $\xi$ contributes only a smoothing symbol. Since each $q_{-m-j}$ can be deformed by the principal homotopy while $r$ remains on $\Gamma_\theta$, for large $\|\xi\|$ we have exactly
\begin{equation}
 b-\sum_{j=0}^{L-1}\chi_0 b_{mz-j}
 =
 \frac{1}{2\pi i}\int_{\Gamma_\theta}
 \lambda^zR_L(x,\xi,\lambda)\,d\lambda.
 \label{pseudo:eq-resto-potencia-contorno-fijo}
\end{equation}

Write $X=\langle\xi\rangle^m$ and let $r_*>0$ be the radius of the fixed curved part. If $K\Subset\{\operatorname{Re}z<0\}$, then after $\partial_z^h$, \eqref{pseudo:eq-resto-resolvente-contorno-fijo} and splitting each branch into $r_*\leq|\lambda|\leq X$ and $|\lambda|\geq X$ give
\begin{align}
 &\left\|D_x^\alpha D_\xi^\beta\partial_z^h
 \left(b-\sum_{j=0}^{L-1}\chi_0 b_{mz-j}\right)\right\| \notag\\
 &\quad\leq
 C_{K,\alpha,\beta,L,h}\langle\xi\rangle^{-L-|\beta|}
 \left[
 X^{-1}\int_{r_*}^{X}
 r^{\operatorname{Re}z}(1+|\log r|)^h\,dr
 {}+{}
 \int_X^\infty
 r^{\operatorname{Re}z-1}(1+\log r)^h\,dr
 \right].
 \label{pseudo:eq-integrales-resto-contorno-fijo}
\end{align}
If $X<r_*$, absorb this region into the constant. For $X\geq r_*$, bound the first integral by separating the cases $\operatorname{Re}z>-1$, $=-1$, and $<-1$; integrate the second in $[X,\infty)$. Uniformly for $z\in K$ and every $\delta>0$, the result is
\begin{equation}
 \left\|D_x^\alpha D_\xi^\beta\partial_z^h
 \left(b-\sum_{j=0}^{L-1}\chi_0 b_{mz-j}\right)\right\|
 \leq C_{K,\alpha,\beta,L,h,\delta}
 \left(
 \langle\xi\rangle^{m\operatorname{Re}z-L-|\beta|+\delta}
 +
 \langle\xi\rangle^{-m-L-|\beta|+\delta}
 \right).
 \label{pseudo:eq-dos-ordenes-resto-contorno-fijo}
\end{equation}
Here we used $(1+|\log\langle\xi\rangle|)^h
\leq C_{h,\delta}\langle\xi\rangle^\delta$; this inequality also covers the borderline case $\operatorname{Re}z=-1$.

To prove the exact classical order, first fix $h=0$ and
$N\in\mathbb N_0$. Choose $L$ with
\[
 L\geq N+1,\qquad
 -m-L+\frac{1}{2}
 \leq m\min_{z\in K}\operatorname{Re}z-N.
\]
Apply \eqref{pseudo:eq-dos-ordenes-resto-contorno-fijo} with
$\delta=\displaystyle\frac{1}{2}$. The first exponent does not exceed
$m\operatorname{Re}z-N-|\beta|-\displaystyle\frac{1}{2}$ and the second
does not exceed $m\operatorname{Re}z-N-|\beta|$. Since
$\langle\xi\rangle\geq1$, this gives
\[
 \left\|D_x^\alpha D_\xi^\beta
 \left(b-\sum_{j=0}^{L-1}\chi_0 b_{mz-j}\right)\right\|
 \leq C_{K,\alpha,\beta,L}
       \langle\xi\rangle^{m\operatorname{Re}z-N-|\beta|}.
\]
By homogeneity and uniform bounds on the compact set of unit directions
and parameters $z\in K$, each component, without derivatives with
respect to $z$, satisfies
\[
 \|D_x^\alpha D_\xi^\beta(\chi_0 b_{mz-j})\|
 \leq C_{K,\alpha,\beta,j}
       \langle\xi\rangle^{m\operatorname{Re}z-j-|\beta|}.
\]
Derivatives of the cutoff are supported in a bounded frequency region
and are absorbed into the same constant. Adding the terms with
$N\leq j<L$ to the remainder estimate gives
\[
 \left\|D_x^\alpha D_\xi^\beta
 \left(b-\sum_{j=0}^{N-1}\chi_0 b_{mz-j}\right)\right\|
 \leq C_{K,\alpha,\beta,N}
       \langle\xi\rangle^{m\operatorname{Re}z-N-|\beta|}.
\]
This is the estimate of order $m\operatorname{Re}z-N$ required for a
classical symbol of complex order $mz$.

Now given $N,h,\delta$ and $K$, choose $L\geq N$ so large that
\[
 -m-L+\frac{\delta}{2}
 \leq m\min_{z\in K}\operatorname{Re}z-N.
\]
Apply \eqref{pseudo:eq-dos-ordenes-resto-contorno-fijo} with $\displaystyle\frac{\delta}{2}$. Both terms are bounded by $C\langle\xi\rangle^{m\operatorname{Re}z-N-|\beta|+\delta}$. The finite terms $\displaystyle \sum_{j=N}^{L-1}\chi_0b_{mz-j}$ satisfy the same bound by \eqref{pseudo:eq-holomorfia-componentes-potencias}. The identity \[
 b-\sum_{j=0}^{N-1}\chi_0b_{mz-j}
 =
 \left(b-\sum_{j=0}^{L-1}\chi_0b_{mz-j}\right)
 +\sum_{j=N}^{L-1}\chi_0b_{mz-j}
\] proves \eqref{pseudo:eq-holomorfia-simbolo-completo-potencias}. The rescaling in \eqref{pseudo:eq-componentes-simbolo-potencia} and the majorant \eqref{pseudo:eq-dominante-integral-componentes-potencia} in turn prove \eqref{pseudo:eq-holomorfia-componentes-potencias}. Thus \[
b(x,\xi;z)
\sim
\sum_{j=0}^\infty b_{mz-j}(x,\xi;z)
\] is a classical symbol of complex order $mz$ and is holomorphic at the first two levels of Definition~\ref{pseudo:def-holomorfia-familia-orden-variable}.

The difference between the exact resolvent and the parametrix satisfies, for each $Q$ and each finite number of kernel derivatives, the bound \eqref{pseudo:eq-residuos-resolvente-decaimiento-cuantificado}. Once $\displaystyle Q>\displaystyle\max_{z\in K}\operatorname{Re}z+1$ is chosen, its product with $|\lambda^z|(1+|\log\lambda|)^h$ is integrable on the infinite branches. The integral therefore has a smooth holomorphic kernel with all seminorms uniform on $K$ and does not alter the classical symbol obtained. The symbol estimates, mapping theorem, and differentiation under the integral give the bound \eqref{pseudo:eq-holomorfia-operatorial-potencias}. To verify holomorphy in this norm, if the segment $[z,z+h]$ remains in a compact neighborhood, use Taylor's identity with integral remainder
\[
 b(z+h)-b(z)-h\partial_zb(z)
 =h^2\int_0^1(1-t)\partial_z^2b(z+th)\,dt.
\]
Each symbol seminorm of the right-hand side is at most $C|h|^2$ by the preceding bounds. Quantization is continuous from these seminorms into $\mathcal L(H^r(M,\mathbf{E}),H^{r-\tau_K-\delta}(M,\mathbf{E}))$, so the difference quotient converges there to $\partial_zA^z$. The same identity applied to each seminorm of $\Gamma(\mathbf{E})$ proves holomorphy on $\Gamma(\mathbf{E})$ itself. This establishes all three levels when $\operatorname{Re}z<0$.

For arbitrary $z$, choose $k$ as in \eqref{pseudo:eq-extension-potencias}. The composition calculus gives
\[
A^kA^{z-k}
\in
\Psi_{\mathrm{cl}}^{mk+m(z-k)}
=
\Psi_{\mathrm{cl}}^{mz}.
\]
For a compact set $K\Subset\mathbb C$, choose a single integer $k$ such that $\operatorname{Re}(z-k)\leq-\eta<0$ for every $z\in K$. Composition with the fixed operator $A^k$ is continuous in the symbol topologies and the operator spaces of \eqref{pseudo:eq-holomorfia-operatorial-potencias}; hence it transfers all uniform bounds already proved to $A^z$. The group law proves independence of $k$. If $w=z-k$ has negative real part, taking $j=0$ in \eqref{pseudo:eq-componentes-simbolo-potencia}, the holomorphic functional calculus of the matrix $a_m(x,\xi)$ gives
\[
b_{mw}(x,\xi;w)
=
\frac{1}{2\pi i}
\int_{\Gamma_{\theta,\xi}}
\lambda^w(\lambda I-a_m(x,\xi))^{-1}\,d\lambda
=a_m(x,\xi)^w.
\]
Multiplicativity of the principal symbol under composition then yields $\boldsymbol{\sigma}_{mz}^\Psi(A^z)=a_m^ka_m^w=a_m^z$ for every $z$. The mapping property is that of a pseudodifferential operator of real order $m\operatorname{Re}z$.
\end{proof}

\begin{proposition}[Closed realizations of powers]
\label{pseudo:prop-dominios-potencias-complejas}
Let $(M,\mathbf{g})$ be a smooth closed Riemannian manifold; let $\mathbf{E}\longrightarrow M$ be a complex vector bundle of finite rank equipped with a Hermitian bundle metric; let $A\in\Psi^m_{\mathrm{cl}}(M;\mathbf{E},\mathbf{E})$ be elliptic and invertible, with $m>0$ and a fixed admissible spectral cut; and let $z\in\mathbb C$. Let $d=m\operatorname{Re}z$. If $d>0$, the closure in $L^2(M,\mathbf{E})$ of $A^z|_{\Gamma(M,\mathbf{E})}$ has domain $H^d(M,\mathbf{E})$, and there exist $c_z,C_z>0$ such that
\begin{equation}
 c_z\|\mathbf{u}\|_{H^d(M,\mathbf{E})}
 \leq \|\mathbf{u}\|_{L^2(M,\mathbf{E})}+\|A^z\mathbf{u}\|_{L^2(M,\mathbf{E})}
 \leq C_z\|\mathbf{u}\|_{H^d(M,\mathbf{E})},
 \qquad \mathbf{u}\in\Gamma(M,\mathbf{E}).
 \label{pseudo:eq-equivalencia-grafica-potencia}
\end{equation}
This realization is closed, and $\Gamma(M,\mathbf{E})$ is a core. If $d\leq0$, $A^z|_{\Gamma(M,\mathbf{E})}$ extends uniquely to a bounded operator on $L^2(M,\mathbf{E})$; this is its realization, its domain is all of $L^2(M,\mathbf{E})$, and $\Gamma(M,\mathbf{E})$ is also a core.
\end{proposition}

\begin{proof}
Suppose $d>0$. The preceding theorem gives continuous maps $A^z:H^d(M,\mathbf{E})\to L^2(M,\mathbf{E})$ and $A^{-z}:L^2(M,\mathbf{E})\to H^d(M,\mathbf{E})$. The group law, first on $\Gamma(\mathbf{E})$, provides the exact identity $\mathbf{u}=A^{-z}A^z\mathbf{u}$. Thus \[
\begin{aligned}
\|\mathbf{u}\|_{H^d(M,\mathbf{E})}
&\leq\|A^{-z}\|_{\mathcal L(L^2(M,\mathbf{E}),H^d(M,\mathbf{E}))}
\|A^z\mathbf{u}\|_{L^2(M,\mathbf{E})},\\
\|A^z\mathbf{u}\|_{L^2(M,\mathbf{E})}
&\leq\|A^z\|_{\mathcal L(H^d(M,\mathbf{E}),L^2(M,\mathbf{E}))}
\|\mathbf{u}\|_{H^d(M,\mathbf{E})}.
\end{aligned}
\], giving \eqref{pseudo:eq-equivalencia-grafica-potencia}. If $\mathbf{u}_j\to \mathbf{u}$ and $A^z\mathbf{u}_j\to \mathbf{f}$ in $L^2(M,\mathbf{E})$, the first inequality applied to $\mathbf{u}_j-\mathbf{u}_k$ proves that $(\mathbf{u}_j)$ converges in $H^d(M,\mathbf{E})$; continuity $A^z:H^d(M,\mathbf{E})\to L^2(M,\mathbf{E})$ gives $A^z\mathbf{u}=\mathbf{f}$. Thus the realization with domain $H^d(M,\mathbf{E})$ is closed. For $\mathbf{u}\in H^d(M,\mathbf{E})$, a sequence of smooth sections converging to $\mathbf{u}$ in $H^d(M,\mathbf{E})$ also converges in graph norm by the second inequality. Consequently, the closure of the restriction to $\Gamma(\mathbf{E})$ is precisely this realization, and $\Gamma(\mathbf{E})$ is a core.

If $d\leq0$, the mapping property gives $A^z:L^2(M,\mathbf{E})\to H^{-d}(M,\mathbf{E})\hookrightarrow L^2(M,\mathbf{E})$ when $d<0$, and $A^z:L^2(M,\mathbf{E})\to L^2(M,\mathbf{E})$ when $d=0$. Density of $\Gamma(\mathbf{E})$ in $L^2(M,\mathbf{E})$ proves uniqueness of the extension. If $\mathbf{u}_j\in\Gamma(\mathbf{E})$ converges to $\mathbf{u}$ in $L^2(M,\mathbf{E})$, continuity implies $A^z\mathbf{u}_j\to A^z\mathbf{u}$ in $L^2(M,\mathbf{E})$, so smooth sections are dense also in graph norm.
\end{proof}

This proof follows the parameter-dependent classes and sums, the resolvent, and the holomorphic families presented in \cite[Chapter~II, §§9.1, 9.4 and 11.1--11.4]{Shubin2001}. Historically, the result and the method of integrating the resolvent parametrix are due to Seeley \cite{Seeley1967}.

\subsection{Agreement with the spectral calculus}

The contour construction also applies to non-self-adjoint operators. When $A$ is self-adjoint and positive, each spectral projection reduces the contour to the scalar Cauchy formula. This identifies Seeley's powers with the spectral calculus used in Chapter~\ref{cap:sobolev-fraccionario-haces-nucleo-calor}.

\begin{proposition}[Contour calculus, spectral calculus, and the semigroup]
\label{pseudo:prop-coincidencia-calculos}
Let $(M,\mathbf{g})$ be a smooth closed Riemannian manifold; let $\mathbf{E}\longrightarrow M$ be a complex vector bundle of finite rank equipped with a Hermitian bundle metric; and let $A\in\Psi^m_{\mathrm{cl}}(M;\mathbf{E},\mathbf{E})$ be elliptic and self-adjoint on $L^2(M,\mathbf{E})$, with $m>0$ and $A\geq cI$ for some $c>0$. Take the negative real axis as the cut. Then the powers in Definition~\ref{pseudo:def-potencias-complejas} agree with the powers from the spectral theorem. As closed operators,
\begin{equation}
A^z
=
\int_{[c,\infty)}\lambda^z\,dE_A(\lambda).
\label{pseudo:eq-potencia-espectral}
\end{equation}
Moreover, for $\operatorname{Re}s>0$,
\begin{equation}
A^{-s}
=
\frac{1}{\Gamma(s)}
\int_0^\infty t^{s-1}e^{-tA}\,dt,
\label{pseudo:eq-mellin-potencias-calor}
\end{equation}
where the Bochner integral converges in $\mathcal L(L^2(M,\mathbf{E}),L^2(M,\mathbf{E}))$.
\end{proposition}

\begin{proof}

First justify the cut from the hypotheses. The adjoint formula makes
the principal Fourier symbol $a_m(x,\xi)$ Hermitian. Fix a chart with an
orthonormal frame, a nonzero constant covector $\xi$, a vector $v$, and a
real $\chi\in C_c^\infty$ in that chart. Extend
$\mathbf u_t(x)=\chi(x)e^{itx\cdot\xi}v$ by zero to $M$. Choose a cutoff
$\zeta$ equal to one on a neighborhood of $\operatorname{supp}\chi$.
In coordinates, $\zeta A\zeta=\operatorname{Op}(a)+R$, where $R$ has a
smooth kernel and the principal component is $a_m$ on that support.
The following integral converges absolutely because $\widehat\chi$ is
a Schwartz function:
\[
 e^{-itx\cdot\xi}\operatorname{Op}(a)
       (\chi e^{it(\cdot)\cdot\xi}v)(x)
 =(2\pi)^{-\frac{n}{2}}
   \int_{\mathbb R^n}e^{ix\cdot\eta}
          a(x,t\xi+\eta)\widehat\chi(\eta)v\,d\eta.
\]
Set $r=\|\xi\|>0$ and take $t$ so large that
$\displaystyle\frac{tr}{2}$ exceeds the frequency beyond which the
classical expansion is used. When
$\|\eta\|\leq\displaystyle\frac{tr}{2}$, the segment
$\xi+\theta\displaystyle\frac{\eta}{t}$, $0\leq\theta\leq1$, stays at
distance at least $\displaystyle\frac{r}{2}$ from zero. Homogeneity of
$a_m$, the fundamental theorem of calculus on this segment, and the
estimate of order $m-1$ for $a-a_m$ give, uniformly for
$x\in\operatorname{supp}\chi$,
\[
 \|t^{-m}a(x,t\xi+\eta)-a_m(x,\xi)\|
 \leq\frac{C}{t}(1+\|\eta\|).
\]
On the region $\|\eta\|>\displaystyle\frac{tr}{2}$, the symbol estimate
of order $m$, with $t\geq1$ and $m>0$, gives
\[
 \|t^{-m}a(x,t\xi+\eta)-a_m(x,\xi)\|
 \leq C\langle\eta\rangle^m,
 \qquad
 1\leq\frac{2\|\eta\|}{tr}.
\]
Thus the integral on the first region is bounded by
$\displaystyle\frac{C}{t}\int(1+\|\eta\|)|\widehat\chi(\eta)|\,d\eta$
and the second by
$\displaystyle\frac{C}{t}\int\langle\eta\rangle^{m+1}|\widehat\chi(\eta)|\,d\eta$.
With
\[
 I_\chi:=\int_{\mathbb R^n}
             \langle\eta\rangle^{m+1}|\widehat\chi(\eta)|\,d\eta<\infty,
\]
Fourier inversion yields
\[
 \sup_{x\in\operatorname{supp}\chi}
 \left\|t^{-m}e^{-itx\cdot\xi}\operatorname{Op}(a)
          (\chi e^{it(\cdot)\cdot\xi}v)(x)
          -a_m(x,\xi)\chi(x)v\right\|
 \leq\frac{CI_\chi\|v\|}{t}.
\]
Let $k(x,y)$ be the coordinate kernel of $R$, including the integration
density. One integration by parts in $y$, with no boundary terms because
$\chi$ has compact support, gives
\[
 \sup_{x\in\operatorname{supp}\chi}\|R\mathbf u_t(x)\|
 \leq\frac{\|v\|}{tr^2}
      \sup_{x\in\operatorname{supp}\chi}
      \int_{\mathbb R^n}
        \|(\xi\cdot\partial_y)(k(x,y)\chi(y))\|\,dy
 \leq\frac{C_\chi\|v\|}{t}.
\]
The final constant is finite by smoothness of the kernel on the relevant
compact set. Pairing with $\chi v$ and integrating over its support gives
\[
 \left|t^{-m}(A\mathbf u_t,\mathbf u_t)_{L^2(M,\mathbf E)}
 -\int_M\chi(x)^2\langle a_m(x,\xi)v,v\rangle
                       \,d\lambda_{\mathbf g}(x)\right|
 \leq C_{\chi,v,\xi}\left(\frac{1}{t}+\frac{1}{t^{m+1}}\right)
 \longrightarrow0.
\]
The limit follows from $m>0$ and the independence of the constants from
$t$. Since $A\geq cI$, the limit is nonnegative. A negative value of
$\langle a_m(x,\xi)v,v\rangle$ at a point would persist on a small
neighborhood and contradict the formula with $\chi$ supported there.
Thus $a_m$ is nonnegative; ellipticity makes it positive definite.
Compactness of the unit cosphere gives a uniform positive lower bound.
The spectrum of $A$ lies in $[c,\infty)$, so the negative real axis avoids
both the principal-symbol spectrum and the operator spectrum, and
$0\in\rho(A)$. It is therefore an admissible cut.

Proposition~\ref{pseudo:prop-resolvente-compacto-espectro-discreto} provides an orthonormal basis $(\boldsymbol{\phi}_j)$ of smooth eigensections, with $A\boldsymbol{\phi}_j=\lambda_j\boldsymbol{\phi}_j$, with the finite-basis convention when applicable. For $\operatorname{Re}z<0$, Cauchy's formula applied to each eigenvalue gives
\[
\frac{1}{2\pi i}
\int_{\Gamma_\pi}
\lambda^z(\lambda I-A)^{-1}\boldsymbol{\phi}_j\,d\lambda
=
\left(
\frac{1}{2\pi i}
\int_{\Gamma_\pi}
\frac{\lambda^z}{\lambda-\lambda_j}\,d\lambda
\right)\boldsymbol{\phi}_j
=
\lambda_j^z\boldsymbol{\phi}_j.
\]
For arbitrary $z$, choose $k\in\mathbb N$ with $\operatorname{Re}(z-k)<0$ and use $A^z=A^kA^{z-k}$; since $A^k\boldsymbol{\phi}_j=\lambda_j^k\boldsymbol{\phi}_j$, this gives $A^z\boldsymbol{\phi}_j=\lambda_j^z\boldsymbol{\phi}_j$. By linearity, both constructions agree on the space $\mathcal F$ of finite spectral sums.

The spectral power has domain
\begin{equation}
 \mathcal D(A^z_{\mathrm{esp}})
 =\left\{\mathbf{u}=\sum_{j=1}^{\infty} u_j\boldsymbol{\phi}_j\in L^2(M,\mathbf{E})
 \ \middle|\
 \sum_{j=1}^{\infty}\lambda_j^{2\operatorname{Re}z}|u_j|^2<\infty\right\}.
 \label{pseudo:eq-dominio-potencia-espectral}
\end{equation}
If $\mathbf{u}$ belongs to this domain and \[
 \mathbf{u}_N:=\sum_{\lambda_j\leq N}u_j\boldsymbol{\phi}_j,
\], then $\mathbf{u}_N\in\mathcal F\subseteq\Gamma(M,\mathbf{E})$ and
\[
 \|\mathbf{u}-\mathbf{u}_N\|_{L^2(M,\mathbf{E})}^2
 +\|A^z_{\mathrm{esp}}(\mathbf{u}-\mathbf{u}_N)\|_{L^2(M,\mathbf{E})}^2
=
\sum_{\lambda_j>N}(1+|\lambda_j^z|^2)|u_j|^2
\longrightarrow0.
\]
Thus $\mathcal F$ is a core for the spectral power.

Let us prove that it is also a core for the pseudodifferential power. If $d=m\operatorname{Re}z\leq0$, this power is bounded on $L^2(M,\mathbf{E})$ by Proposition~\ref{pseudo:prop-dominios-potencias-complejas}; density of $\mathcal F$ in $L^2(M,\mathbf{E})$ then implies density in graph norm. Suppose $d>0$ and choose
\[
N_0:=\max\left\{1,\left\lceil\frac{d}{m}\right\rceil\right\},
\qquad mN_0\geq d.
\]
If $\displaystyle \mathbf{v}=\displaystyle\sum_{j=1}^{\infty}v_j\boldsymbol{\phi}_j\in\Gamma(M,\mathbf{E})$ and $\displaystyle \mathbf{v}_N=\displaystyle\sum_{\lambda_j\leq N}v_j\boldsymbol{\phi}_j$, then $A^{N_0}\mathbf{v}\in L^2(M,\mathbf{E})$. Since $\mathbf{v}$ belongs to the domain of $A^{N_0}$ and $A$ is self-adjoint, the coefficient of $A^{N_0}\mathbf{v}$ in $\boldsymbol{\phi}_j$ is
\[
 (A^{N_0}\mathbf{v},\boldsymbol{\phi}_j)_{L^2(M,\mathbf{E})}
 =(\mathbf{v},A^{N_0}\boldsymbol{\phi}_j)_{L^2(M,\mathbf{E})}=\lambda_j^{N_0}v_j.
\]
By Parseval,
\[
 \|\mathbf{v}-\mathbf{v}_N\|_{L^2(M,\mathbf{E})}^2+
 \|A^{N_0}(\mathbf{v}-\mathbf{v}_N)\|_{L^2(M,\mathbf{E})}^2
 =\sum_{\lambda_j>N}(1+\lambda_j^{2N_0})|v_j|^2
 \longrightarrow0.
\]
Graph-norm equivalence \eqref{pseudo:eq-equivalencia-grafica-potencia}, applied to exponent $N_0$, gives $\mathbf{v}_N\to \mathbf{v}$ in $H^{mN_0}(M,\mathbf{E})$ and, by the continuous inclusion, in $H^d(M,\mathbf{E})$.

Now let $\mathbf{u}\in H^d(M,\mathbf{E})$ and $\varepsilon>0$. By density of $\Gamma(M,\mathbf{E})$ in $H^d(M,\mathbf{E})$, choose $\mathbf{v}\in\Gamma(M,\mathbf{E})$ such that
\[
\|\mathbf{u}-\mathbf{v}\|_{H^d(M,\mathbf{E})}<\frac{\varepsilon}{2}.
\]
For this $\mathbf{v}$, the convergence just proved allows $N\in\mathbb N$ to be chosen so that
\[
\|\mathbf{v}-\mathbf{v}_N\|_{H^d(M,\mathbf{E})}<\frac{\varepsilon}{2}.
\]
Since $\mathbf{v}_N\in\mathcal F$, the triangle inequality gives
\[
\|\mathbf{u}-\mathbf{v}_N\|_{H^d(M,\mathbf{E})}<\varepsilon.
\]
Thus $\mathcal F$ is dense in $H^d(M,\mathbf{E})$. Equivalence \eqref{pseudo:eq-equivalencia-grafica-potencia} for $z$ shows that this is also density in the graph norm of the pseudodifferential power. Hence $\mathcal F$ is a core for both closed realizations. Since they agree on $\mathcal F$, their closures coincide, proving \eqref{pseudo:eq-potencia-espectral} without presupposing equality on all of $\Gamma(\mathbf{E})$.

The scalar identity
\[
\lambda^{-s}
=
\frac{1}{\Gamma(s)}
\int_0^\infty t^{s-1}e^{-t\lambda}\,dt,
\qquad \lambda\geq c,
\]
is transferred by the spectral theorem. Near zero, $t^{\operatorname{Re}s-1}$ is a majorant, and at infinity $t^{\operatorname{Re}s-1}e^{-ct}$ is a majorant. More precisely,
\[
 \|t^{s-1}e^{-tA}\|_{\mathcal L(L^2(M,\mathbf{E}))}
 \leq t^{\operatorname{Re}s-1}
 \quad(0<t\leq1),
 \qquad
 \|t^{s-1}e^{-tA}\|_{\mathcal L(L^2(M,\mathbf{E}))}
 \leq t^{\operatorname{Re}s-1}e^{-ct}
 \quad(t\geq1).
\]
The first function is integrable on $(0,1]$ because $\operatorname{Re}s>0$, and the second on $[1,\infty)$. The Bochner integral therefore converges in $\mathcal L(L^2(M,\mathbf{E}))$. Applying both expressions to each spectral projection and using the same majorants allows the integral and spectral calculus to be interchanged, proving \eqref{pseudo:eq-mellin-potencias-calor}.
\end{proof}

\subsection{Bessel potentials, the fractional Laplacian, and heat}

The preceding identification returns us to the operators that motivated this chapter. Powers of the connection Laplacian, previously defined by the spectral theorem and the semigroup, now have classical symbols and can be studied within the same microlocal calculus.

Now let $(M,\mathbf{g})$ be a smooth closed Riemannian manifold, and let $\mathbf{E}\longrightarrow M$ be a complex vector bundle of finite rank equipped with a Hermitian bundle metric $\mathbf{h}_{\mathbf{E}}$ and a compatible connection $\nabla^{\mathbf{E}}$. Denote its nonnegative self-adjoint realization by $L_{\mathbf{E}}=(\nabla^{\mathbf{E}})_h^*\nabla^{\mathbf{E}}$. Although the intrinsic differential symbol of $L_{\mathbf{E}}$ is $-\|\xi\|_{\mathbf{g}}^2I$, its principal Fourier symbol is $\|\xi\|_{\mathbf{g}}^2I$. Thus the operator \[
B_{\mathbf{E}}:=I+L_{\mathbf{E}}
\] is elliptic, self-adjoint, satisfies $B_{\mathbf{E}}\geq I$, and admits the negative axis as a spectral cut.

\begin{corollary}[Bessel powers and the fractional Laplacian]
\label{pseudo:cor-bessel-laplaciano-fraccionario}
\index{Bessel potential!pseudodifferential description}
\index{fractional Laplacian}
Let $(M,\mathbf{g})$ be a smooth closed Riemannian manifold, and let $\mathbf{E}\longrightarrow M$ be a complex vector bundle of finite rank equipped with a Hermitian bundle metric $\mathbf{h}_{\mathbf{E}}$ and a compatible connection $\nabla^{\mathbf{E}}$. Let $L_{\mathbf{E}}=(\nabla^{\mathbf{E}})_h^*\nabla^{\mathbf{E}}$ and $B_{\mathbf{E}}=I+L_{\mathbf{E}}$. For every $z\in\mathbb C$,
\[
B_{\mathbf{E}}^z\in\Psi^{2z}_{\mathrm{cl}}(M;\mathbf{E},\mathbf{E}),
\qquad
\boldsymbol{\sigma}_{2z}^\Psi(B_{\mathbf{E}}^z)(x,\xi)
=\|\xi\|_{\mathbf{g}}^{2z}I_{\mathbf{E}_x}.
\]
If $s>0$, the Bessel potential constructed through the semigroup satisfies
\begin{equation}
J_{s,2}^{\mathbf{E}}
=(I+L_{\mathbf{E}})^{-\frac{s}{2}}
=B_{\mathbf{E}}^{-\frac{s}{2}}
\in\Psi^{-s}_{\mathrm{cl}}(M;\mathbf{E},\mathbf{E}).
\label{pseudo:eq-bessel-pseudo}
\end{equation}
In particular, for every $r\in\mathbb R$,
\[
B_{\mathbf{E}}^{-\frac{s}{2}}:H^r(M,\mathbf{E})\longrightarrow H^{r+s}(M,\mathbf{E})
\]
is an isomorphism with inverse $B_{\mathbf{E}}^{\displaystyle\frac{s}{2}}$. More precisely, for every $s\in\mathbb R$ there exist constants $c_s,C_s>0$ such that, for $\mathbf{u}\in\Gamma(M,\mathbf{E})$,
\begin{equation}
c_s\|\mathbf{u}\|_{H^s(M,\mathbf{E})}
\leq
\|B_{\mathbf{E}}^{\frac{s}{2}}\mathbf{u}\|_{L^2(M,\mathbf{E})}
\leq
C_s\|\mathbf{u}\|_{H^s(M,\mathbf{E})}.
\label{pseudo:eq-equivalencia-norma-bessel}
\end{equation}
The constants depend on $s$, the metric $\mathbf{g}$, the Hermitian metric and connection of $\mathbf{E}$, and the finite cover, trivializations, and partition of unity chosen for the localized Sobolev norm; once these data are fixed, they are uniform in $\mathbf{u}$.

Let $\Pi_0$ be the orthogonal projection onto $\ker L_{\mathbf{E}}$. For $s>0$,
\begin{equation}
L_{\mathbf{E}}^{\frac{s}{2}}
=(L_{\mathbf{E}}+\Pi_0)^{\frac{s}{2}}-\Pi_0
\in\Psi^s_{\mathrm{cl}}(M;\mathbf{E},\mathbf{E}),
\label{pseudo:eq-laplaciano-fraccionario}
\end{equation}
and the reduced negative power on $(\ker L_{\mathbf{E}})^\perp$ is
\begin{equation}
L_{\mathbf{E},\perp}^{-\frac{s}{2}}
=(L_{\mathbf{E}}+\Pi_0)^{-\frac{s}{2}}-\Pi_0
\in\Psi^{-s}_{\mathrm{cl}}(M;\mathbf{E},\mathbf{E}).
\label{pseudo:eq-laplaciano-potencia-reducida}
\end{equation}
Both have principal symbol $\|\xi\|_{\mathbf{g}}^{\pm s}I_{\mathbf{E}_x}$.
\end{corollary}

\begin{proof}
Theorem~\ref{pseudo:teo-potencias-complejas-clasicas} applied to $B_{\mathbf{E}}$ gives the first assertion. Proposition~\ref{pseudo:prop-coincidencia-calculos} and the spectral characterization in Theorem~\ref{teo:caracterizacion-espectral-bessel-haces} prove \eqref{pseudo:eq-bessel-pseudo}. The mapping properties of $B_{\mathbf{E}}^{\pm\displaystyle\frac{s}{2}}$ and the group law give the Sobolev isomorphism. The spectral identity
\[
\|B_{\mathbf{E}}^{\frac{s}{2}}\mathbf{u}\|_{L^2(M,\mathbf{E})}^2
=\sum_{j=1}^{\infty}(1+\lambda_j)^s|u_j|^2
\]
is the norm of the spectral scale in the preceding chapter. Equivalence of this scale with the localized scale, applied in both directions, gives exactly \eqref{pseudo:eq-equivalencia-norma-bessel} and the stated dependencies. Theorem~\ref{pseudo:teo-potencias-complejas-clasicas} adds the pseudodifferential description of the same scale.

The kernel of $L_{\mathbf{E}}$ is finite-dimensional and consists of smooth sections. Thus $\Pi_0$ has smooth kernel and belongs to $\Psi^{-\infty}(M;\mathbf{E},\mathbf{E})$. The operator $L_{\mathbf{E}}+\Pi_0$ is strictly positive, invertible, and has the same principal symbol as $L_{\mathbf{E}}$. On $\ker L_{\mathbf{E}}$ its eigenvalue is one, while on the complement it agrees with $L_{\mathbf{E}}$. The spectral calculus then gives identities \eqref{pseudo:eq-laplaciano-fraccionario} and \eqref{pseudo:eq-laplaciano-potencia-reducida}. Subtracting $\Pi_0$ does not change the principal symbol.
\end{proof}

On a manifold without boundary and of bounded geometry, Theorems~\ref{teo:yoshida-bessel-covariante-geometria-acotada} and~\ref{teo:triebel-bessel-localizacion-geometria-acotada} already identify, under their stated hypotheses, the spectral scale of $I+L_{\mathbf{E}}$ with the scale defined by uniform localization. Both inequalities in \eqref{pseudo:eq-equivalencia-norma-bessel} then retain the same form, with constants depending on $s$ and quantitative bounded-geometry bounds, but not on the chart. Theorem~\ref{pseudo:teo-mapeo-uniforme-sobolev} explains the mapping properties of the resulting uniform operators. We do not, however, assert that every complex power automatically belongs to $\Psi^{2z}_{\mathrm u}(M;\mathbf{E},\mathbf{E})$: that conclusion requires a uniform parameter-dependent calculus and global resolvent estimates at infinity, hypotheses that do not follow from the local construction. On a closed manifold, finiteness of the atlas and compactness of the resolvent remove precisely this difficulty and yield the full result proved above.

For the trivial bundle over $\mathbb R^n$ and $s>0$, the same identities can be read directly in Fourier variables. With the unitary normalization fixed in this chapter, the fractional Laplacian is defined exactly as the multiplier
\[
\widehat{(-\Delta)^{\frac{s}{2}}u}(\xi)
=\|\xi\|^s\widehat u(\xi),
\qquad u\in\mathcal S(\mathbb R^n),
\]
with closed domain given by
\[
\left\{
u\in L^2(\mathbb R^n)
\;\middle|\;
\|\xi\|^s\widehat u(\xi)\in L^2(\mathbb R^n)
\right\}.
\]
In contrast, \[
(I-\Delta)^{\frac{s}{2}}
=\operatorname{Op}\bigl((1+\|\xi\|^2)^{\frac{s}{2}}\bigr)
\] belongs directly to the global inhomogeneous calculus. To express the pseudodifferential part of the first operator, choose $\chi\in C^\infty(\mathbb R^n)$ vanishing near zero and equal to one for large $\|\xi\|$, and reserve quantization notation for the smooth symbol
\[
\operatorname{Op}\bigl(\chi(\xi)\|\xi\|^s\bigr)
\in\Psi^s_{\mathrm{cl,loc}}(\mathbb R^n;\mathbb C,\mathbb C).
\]
The symbol is classical, but the frequency cutoff does not imply proper support of the operator. The subscript $\mathrm{loc}$ is therefore retained. If a properly supported representative is required, multiply its kernel by a properly supported cutoff equal to one near the diagonal. The difference has smooth kernel because the original kernel is smooth off the diagonal; the resulting representative is not generally the exact Fourier multiplier. The difference between the exact fractional Laplacian and its high-frequency part is supported in a compact frequency set and maps $L^2(\mathbb R^n)$ continuously into $H^N(\mathbb R^n)$ for every $N$, without identifying $\|\xi\|^s$ with a symbol smooth at the origin. On a closed manifold, the frequency-zero difficulty is concentrated in the finite-dimensional space $\ker L_{\mathbf{E}}$, and the formulas with $\Pi_0$ isolate it exactly. At high frequencies, however, the principal symbols of $L_{\mathbf{E}}^{\displaystyle\frac{s}{2}}$ and $(I+L_{\mathbf{E}})^{\displaystyle\frac{s}{2}}$ coincide. Except for exponents producing a differential operator, these powers are nonlocal: their Schwartz kernels are not supported on the diagonal. For $0<s<2$, this nonlocality is also visible in the representation, valid for $\mathbf{u}\in\Gamma(\mathbf{E})$,
\[
L_{\mathbf{E}}^{\frac{s}{2}}\mathbf{u}
=
\frac{s}{2\Gamma(1-\frac{s}{2})}
\int_0^\infty
\frac{\mathbf{u}-e^{-tL_{\mathbf{E}}}\mathbf{u}}{t^{1+\frac{s}{2}}}\,dt.
\]
The integral converges as a Bochner integral in $L^2(M,\mathbf{E})$. For $\mathbf{u}\in\mathcal D(L_{\mathbf{E}})$, the semigroup equation and the fundamental theorem of calculus in $L^2(M,\mathbf{E})$ give the exact identity
\begin{equation}
 \mathbf{u}-e^{-tL_{\mathbf{E}}}\mathbf{u}
 =\int_0^tL_{\mathbf{E}}e^{-\tau L_{\mathbf{E}}}\mathbf{u}\,d\tau.
 \label{pseudo:eq-identidad-incremento-semigrupo}
\end{equation}
Since the semigroup is contractive and commutes with $L_{\mathbf{E}}$ on $\mathcal D(L_{\mathbf{E}})$,
\begin{equation}
 \|\mathbf{u}-e^{-tL_{\mathbf{E}}}\mathbf{u}\|_{L^2(M,\mathbf{E})}
 \leq\int_0^t\|e^{-\tau L_{\mathbf{E}}}L_{\mathbf{E}}\mathbf{u}\|_{L^2(M,\mathbf{E})}\,d\tau
 \leq t\|L_{\mathbf{E}}\mathbf{u}\|_{L^2(M,\mathbf{E})}.
 \label{pseudo:eq-cota-incremento-semigrupo}
\end{equation}
Thus, for $0<t\leq1$, the norm of the integrand is at most $\|L_{\mathbf{E}}\mathbf{u}\|_{L^2(M,\mathbf{E})}t^{-\displaystyle\frac{s}{2}}$, which is integrable because $s<2$. For $t\geq1$, contractivity gives
\[
 \left\|\frac{\mathbf{u}-e^{-tL_{\mathbf{E}}}\mathbf{u}}{t^{1+\frac{s}{2}}}\right\|_{L^2(M,\mathbf{E})}
 \leq2\|\mathbf{u}\|_{L^2(M,\mathbf{E})}t^{-1-\frac{s}{2}},
\]
and this majorant is integrable because $s>0$. The $\ker L_{\mathbf{E}}$ component vanishes exactly in the numerator.

Negative powers and heat are two presentations of the same functional calculus. The first integrates the resolvent against $\lambda^{-s}$; the second integrates the spectrum against $e^{-t\lambda}$. The Mellin transform in \eqref{pseudo:eq-mellin-potencias-calor} identifies them.

\begin{proposition}[Heat smoothing and trace]
\label{pseudo:prop-calor-regularizacion}
\index{heat operator!smoothing}
Let $(M,\mathbf{g})$ be a smooth closed Riemannian manifold, and let $\mathbf{E}\longrightarrow M$ be a complex vector bundle of finite rank equipped with a Hermitian bundle metric $\mathbf{h}_{\mathbf{E}}$ and a compatible connection $\nabla^{\mathbf{E}}$. Let $A\geq0$ be a classical elliptic self-adjoint pseudodifferential operator of order $m>0$ on $\mathbf{E}$. For each $t>0$,
\[
e^{-tA}\in\Psi^{-\infty}(M;\mathbf{E},\mathbf{E}).
\]
Its kernel $\mathbf{K}_A(t,x,y)$ is smooth on $(0,\infty)\times M\times M$. Moreover, $e^{-tA}$ is trace class and
\begin{equation}
\operatorname{Tr}(e^{-tA})
=
\int_{M_x}\operatorname{tr}_{\mathbf{E}_x}\mathbf{K}_A(t,x,x)\,d\lambda_{\mathbf{g}}(x).
\label{pseudo:eq-traza-calor-diagonal}
\end{equation}
\end{proposition}

\begin{proof}
The spectral calculus gives, for each $N\in\mathbb N$,
\[
\|A^Ne^{-tA}\|_{\mathcal L(L^2(M,\mathbf{E}),L^2(M,\mathbf{E}))}
\leq
\sup_{\lambda\geq0}\lambda^Ne^{-t\lambda}<\infty.
\]
Thus $e^{-tA}$ maps $L^2(M,\mathbf{E})$ into $\mathcal D(A^N)$. Elliptic regularity identifies this domain, locally and with equivalent norms, with $H^{mN}(M,\mathbf{E})$. Since $N$ is arbitrary, the operator gains arbitrarily many derivatives. The same argument for the adjoint shows that its kernel is smooth; hence it belongs to $\Psi^{-\infty}(M;\mathbf{E},\mathbf{E})$. For $t$ in a compact subset of $(0,\infty)$, the identities $\partial_t^k e^{-tA}=(-A)^ke^{-tA}$ and the same estimates, with arbitrary $N$, give joint smoothness in $(t,x,y)$. On a compact manifold, a smooth kernel defines a trace-class operator and the diagonal formula \eqref{pseudo:eq-traza-calor-diagonal} holds.
\end{proof}

For the McKean--Singer argument, smoothing, smoothness of the kernel, and the diagonal trace formula suffice; the fine heat asymptotics will be needed only later, in studying the local index density. For $B_{\mathbf{E}}=I+L_{\mathbf{E}}$, we obtain
\[
e^{-tB_{\mathbf{E}}}=e^{-t}e^{-tL_{\mathbf{E}}}
\]
and
\[
B_{\mathbf{E}}^{-\frac{s}{2}}
=
\frac{1}{\Gamma(\frac{s}{2})}
\int_0^\infty
t^{\frac{s}{2}-1}e^{-t}e^{-tL_{\mathbf{E}}}\,dt,
\]
which is exactly the Bochner formula used to define Bessel potentials.

\subsection{From heat to the index}

The chapter's last step shows how the preceding theory enters the index problem. The parametrix gives the Fredholm property, and heat turns the difference between kernel and cokernel into a difference of traces. Let $(M,\mathbf{g})$ be a smooth closed Riemannian manifold; let $\mathbf{E}^+,\mathbf{E}^-\longrightarrow M$ be complex vector bundles of finite rank equipped with Hermitian bundle metrics and Hermitian connections; and let \[
P\in\Psi^m_{\mathrm{cl}}(M;\mathbf{E}^+,\mathbf{E}^-),
\qquad
P\colon\Gamma(M,\mathbf{E}^+)\longrightarrow\Gamma(M,\mathbf{E}^-),
\qquad m>0,
\] be a classical elliptic operator. Its parametrix implies that its closed realization \[
P\colon H^m(M,\mathbf{E}^+)\subset L^2(M,\mathbf{E}^+)\longrightarrow L^2(M,\mathbf{E}^-)
\] is Fredholm. Its Hilbert space adjoint $P^\dagger$ has domain $H^m(M,\mathbf{E}^-)$ and agrees with the closed realization of the formal adjoint $P_h^*$. Thus
\[
\operatorname{ind}P
=\dim\ker P-\dim\ker P^\dagger.
\]
This functional-analytic part of the index is a consequence of Atkinson's theorem and compactness of remainders on a closed manifold; see also \cite[Chapters~3 and~9]{BleeckerBoossIndex}.

\begin{proposition}[McKean--Singer formula]
\label{pseudo:prop-mckean-singer}
\index{McKean--Singer formula@McKean--Singer formula}
Let $(M,\mathbf{g})$ be a smooth closed Riemannian manifold; let $\mathbf{E}^+,\mathbf{E}^-\longrightarrow M$ be complex vector bundles of finite rank equipped with Hermitian bundle metrics; and let $P\in\Psi^m_{\mathrm{cl}}(M;\mathbf{E}^+,\mathbf{E}^-)$ be elliptic, with $m>0$. For every $t>0$,
\begin{equation}
\operatorname{ind}P
=
\operatorname{Tr}(e^{-tP^\dagger P})
-
\operatorname{Tr}(e^{-tPP^\dagger}).
\label{pseudo:eq-mckean-singer}
\end{equation}
\end{proposition}

The spectral proof corresponds to \cite[Lemma~1.6.5]{Gilkey1984}.

\begin{proof}
Write \(H_\pm:=L^2(M,\mathbf E^\pm)\).
The realization of \(P\) is closed with domain \(H^m(M,\mathbf E^+)\),
and its adjoint has domain \(H^m(M,\mathbf E^-)\), as stated before
the proposition. The form
\[
 q_+(u,v):=\langle Pu,Pv\rangle_{H_-},\qquad
 \mathcal D(q_+):=H^m(M,\mathbf E^+),
\]
is nonnegative, densely defined, and closed: its form norm is the
graph norm of \(P\), equivalent to \(H^m\).
Theorem~\ref{teo:representacion-formas-cerradas} gives a nonnegative
self-adjoint operator. The identity
\(q_+(u,v)=\langle f,v\rangle_{H_+}\) for all \(v\in\mathcal D(P)\)
is equivalent, by the definition of the adjoint, to
\(Pu\in\mathcal D(P^\dagger)\) and \(P^\dagger Pu=f\).
Thus this operator is exactly \(P^\dagger P\), with domain
\[
 \mathcal D(P^\dagger P)
 =\{u\in H^m(M,\mathbf E^+)\mid Pu\in H^m(M,\mathbf E^-)\}
 =H^{2m}(M,\mathbf E^+).
\]
The last equality uses elliptic regularity of \(P\) in one direction
and its mapping \(H^{2m}\to H^m\) in the other. Applying the same
argument to \(P^\dagger\), we obtain
\[
 \mathcal D(PP^\dagger)=H^{2m}(M,\mathbf E^-).
\]
On their domains these products agree with the pseudodifferential
expressions \(P_h^*P\) and \(PP_h^*\), whose principal symbols are
positive and invertible for \(\xi\neq0\). They are therefore the
self-adjoint elliptic realizations to which
Proposition~\ref{pseudo:prop-calor-regularizacion} applies; their
semigroups are smoothing and trace class.

Eigenvectors of both products are smooth by elliptic regularity.
If \(P^\dagger Pu=\lambda u\) with \(\lambda>0\), then
\[
 PP^\dagger(Pu)=P(P^\dagger Pu)=\lambda Pu,
 \qquad \|Pu\|_{H_-}^2=\lambda\|u\|_{H_+}^2.
\]
All compositions are defined because \(u\) and \(Pu\) are smooth.
The argument with \(P^\dagger\) gives the inverse
\(\lambda^{-1}P^\dagger\). Positive eigenspaces therefore have
the same dimensions. Both trace series converge absolutely;
grouping them by their common eigenvalues cancels all positive
contributions. Finally,
\[
 \ker(P^\dagger P)=\ker P,\qquad
 \ker(PP^\dagger)=\ker P^\dagger,
\]
since their quadratic forms are, respectively, \(\|Pu\|^2\) and
\(\|P^\dagger v\|^2\). The difference of the zero multiplicities
is \(\operatorname{ind}P\), as asserted.

\end{proof}

The right-hand side of \eqref{pseudo:eq-mckean-singer} is independent of $t$. The diagonal formula \eqref{pseudo:eq-traza-calor-diagonal} interprets both traces and the spectral cancellation; its local limit is expressed through the fine asymptotics of the heat kernel.

The analytic conclusion depends only on the class of the principal symbol. To see this, let us specify which datum remains invariant under the preceding operations. If $\pi:T^*M\setminus0_M\to M$ is the projection and $P$ has order $m$, ellipticity says that \[
\boldsymbol{\sigma}_m^\Psi(P):\pi^*\mathbf{E}^+\longrightarrow\pi^*\mathbf{E}^-
\] is a bundle isomorphism over $T^*M\setminus0_M$. If $R\in\Psi^{m-1}_{\mathrm{cl}}(M;\mathbf{E}^+,\mathbf{E}^-)$, the family $P_t=P+tR$ preserves this principal symbol and is therefore elliptic for every $t\in[0,1]$. Viewed between the appropriate Sobolev spaces, it is a norm-continuous family of Fredholm operators; local constancy of the index on the Fredholm set gives $\operatorname{ind}(P+R)=\operatorname{ind}P$. If $R$ is smoothing, the same conclusion follows directly from Atkinson's theorem, since $R$ is compact. Finally, if \(S^*M=\varnothing\), there are no symbol data away
from the zero section. A closed zero-dimensional manifold is a finite
set, and all operators act between finite-dimensional spaces;
their index is \(\dim L^2(M,\mathbf E^+)-\dim L^2(M,\mathbf E^-)\),
independently of the operator. Now assume \(S^*M\neq\varnothing\)
and fix bundle metrics \(h_+,h_-\). Using the induced norms, set
\[
 M_\sigma:=\max_{(x,\xi)\in S^*M}
             \|\boldsymbol\sigma_m^\Psi(P)(x,\xi)^{-1}\|,
 \qquad \delta:=\frac{1}{2(1+M_\sigma)}.
\]
Compactness makes \(M_\sigma\) finite, also in rank zero.
If another homogeneous principal symbol \(b_m\) satisfies
\(\sup_{S^*M}\|b_m-\boldsymbol\sigma_m^\Psi(P)\|<\delta\), then,
for \(0\leq t\leq1\),
\[
 \|t\boldsymbol\sigma_m^\Psi(P)^{-1}
                      (b_m-\boldsymbol\sigma_m^\Psi(P))\|
 \leq\frac{M_\sigma}{2(1+M_\sigma)}<\frac{1}{2}.
\]
The Neumann series therefore inverts
\[
 \boldsymbol\sigma_m^\Psi(P)+t(b_m-\boldsymbol\sigma_m^\Psi(P))
 =\boldsymbol\sigma_m^\Psi(P)
  [I+t\boldsymbol\sigma_m^\Psi(P)^{-1}
                         (b_m-\boldsymbol\sigma_m^\Psi(P))].
\]
Quantizing the segment gives a norm-continuous family of elliptic
operators with unchanged index.

The analytic problem is thus reduced to the topological class of the principal symbol. The following chapters construct that class and explain how the index is computed from it.

\chapter{Fredholm operators and the analytic index}
\label{cap:fredholm-indice-analitico}

The pseudodifferential calculus of Chapter~\ref{cap:operadores-pseudodiferenciales}
assigns to every elliptic operator on a closed manifold a parametrix whose
remainders are smoothing. Our task now is to isolate the functional-analytic
consequences of this construction and establish the properties of the index
that will be used in passing from the principal symbol to $K$--theory. The
first part is formulated for bounded operators between Banach spaces. The
second returns to Sobolev spaces of sections and shows that the resulting
integer depends on the homotopy class of the elliptic symbol, independently
of the chosen completion.

The functional-analytic exposition follows Chapters~1--3 of
\cite[pp.~3--88]{BleeckerBoossIndex}; the application to elliptic operators
corresponds to Theorem~9.10 and Corollary~9.18 of
\cite[pp.~240--242]{BleeckerBoossIndex}. We give the proofs that make the passage
between the two parts explicit.

\begin{semblanzaHistorica}{Fredholm and the persistence of finite dimension}
Ivar Fredholm discovered that certain infinite-dimensional integral equations retain an alternative governed by finite-dimensional kernels and obstructions. Modern theory abstracts this situation in the notion of a Fredholm operator. Its index does not count all solutions; it measures the stable imbalance between solutions and compatibility conditions. This stability explains why an analytic integer can survive deformations of the operator and ultimately depend only on its principal symbol.
\end{semblanzaHistorica}

\section{Fredholm operators and parametrices modulo compact operators}

Let $X$ and $Y$ be Banach spaces over $\mathbb K$, where
$\mathbb K=\mathbb R$ or $\mathbb C$. Recall
Definition~\ref{def:operador-fredholm-indice}: an operator
$T\in\mathcal L(X,Y)$ is Fredholm when $\ker T$ is finite-dimensional,
$\operatorname{Ran}T$ is closed, and the cokernel
\[
\operatorname{coker}T:=Y/\operatorname{Ran}T
\]
is finite-dimensional. In this case,
\[
\operatorname{ind}T
:=\dim\ker T-\dim\operatorname{coker}T.
\]
We explicitly require the range to be closed because of its analytic role:
it gives the cokernel a Banach space structure and ensures that the inverse
induced by $T$ on a closed complement of its kernel is continuous.
For a bounded operator between Banach spaces, however, closedness also
follows from the finite dimensionality of the algebraic cokernel.

Indeed, suppose that $Y/\operatorname{Ran}T$ is finite-dimensional and
set $K:=\ker T$. The subspace $K$ is closed by the continuity of
$T$, so $X/K$ is Banach by
Proposition~\ref{prop:norma-cociente-banach}. Choose representatives of
a basis of $Y/\operatorname{Ran}T$ and denote their span by $C$.
Then $C$ is finite-dimensional and
$Y=\operatorname{Ran}T\oplus C$ as an algebraic direct sum. Consider
\[
 S\colon (X/K)\times C\longrightarrow Y,
 \qquad S([x],c):=Tx+c.
\]
This operator is well defined because $T$ vanishes on $K$, and the
algebraic direct sum shows that it is linear and bijective. Moreover, for every $k\in K$,
$\|Tx\|_Y=\|T(x-k)\|_Y\leq\|T\|\|x-k\|_X$; taking the infimum yields
\[
 \|S([x],c)\|_Y\leq\|T\|\|[x]\|_{X/K}+\|c\|_Y.
\]
Thus $S$ is continuous when the product is equipped with the sum norm.
Both spaces are Banach, so
Corollary~\ref{teo:inverso-acotado} implies that $S^{-1}$ is continuous.
Consequently,
\[
 \operatorname{Ran}T=S\bigl((X/K)\times\{0\}\bigr)
\]
is closed, since a homeomorphism takes closed sets to closed sets. This
conclusion does not require $\ker T$ to be finite-dimensional.

\begin{lemma}[Complements associated with a Fredholm operator]
\label{lem:descomposiciones-fredholm}
Let $T\in\mathcal L(X,Y)$ be Fredholm. There exist closed subspaces $X_0\subseteq
X$ and $C\subseteq Y$ such that
\[
X=\ker T\oplus X_0,
\qquad
Y=\operatorname{Ran}T\oplus C,
\]
and the restriction
\[
T_0:=T\restriction_{X_0}\colon X_0\longrightarrow\operatorname{Ran}T
\]
is an isomorphism of Banach spaces. Moreover,
$\dim C=\dim\operatorname{coker}T$.
\end{lemma}

\begin{proof}
Every finite-dimensional subspace of a normed space is complemented. Indeed,
if $(e_1,\ldots,e_r)$ is a basis of $\ker T$, the coordinate
functionals on $\ker T$ are continuous, and
Corollary~\ref{cor:extension-hahn-banach-preserva-norma} allows us to extend them to
$\ell_1,\ldots,\ell_r\in X'$. The operator
\[
P_Nx:=\sum_{j=1}^r\ell_j(x)e_j
\]
is a continuous projection onto $\ker T$; therefore,
$X_0:=\ker P_N$ is closed and $X=\ker T\oplus X_0$.

Since $\operatorname{Ran}T$ is closed and has finite codimension, choose
$c_1,\ldots,c_q\in Y$ whose classes form a basis of
$Y/\operatorname{Ran}T$ and set $C:=\operatorname{span}\{c_1,\ldots,c_q\}$.
Every $y\in Y$ has a unique expression $y=r+c$, with
$r\in\operatorname{Ran}T$ and $c\in C$. The projection onto $C$ has closed graph:
if $y_\nu\to y$ and their components satisfy $c_\nu\to c$, then
$y-c\in\operatorname{Ran}T$. Theorem~\ref{teo:grafica-cerrada} shows that the
projection is continuous. Thus $Y=\operatorname{Ran}T\oplus C$ is a topological
direct sum and $\dim C=q=\dim\operatorname{coker}T$.

The restriction $T_0$ is linear, continuous, and injective. It is also surjective:
if $Tx\in\operatorname{Ran}T$ and $x=n+x_0$ is the first decomposition,
then $Tx=Tx_0$. Theorem~\ref{teo:mapeo-abierto-banach} shows that
$T_0^{-1}$ is continuous.
\end{proof}

The central characterization identifies the Fredholm property with
invertibility up to compact errors. This equivalence is Atkinson's theorem.

\begin{theorem}[Atkinson]
\label{teo:atkinson-fredholm}
\index{Atkinson's theorem}
\index{Fredholm operator!Atkinson characterization}
For $T\in\mathcal L(X,Y)$, the following statements are equivalent.
\begin{enumerate}[label=(\roman*)]
\item The operator $T$ is Fredholm.
\item There exists $S\in\mathcal L(Y,X)$ such that
\begin{equation}
ST=I_X-K_X,
\qquad
TS=I_Y-K_Y,
\label{eq:parametrix-atkinson}
\end{equation}
with $K_X\in\mathcal K(X)$ and $K_Y\in\mathcal K(Y)$.
\end{enumerate}
If $X=Y$, the equivalence states that the class of $T$ in the Calkin algebra
\[
\mathcal C(X):=\mathcal L(X)/\mathcal K(X)
\]
is invertible.
\end{theorem}

\begin{proof}
First suppose that $T$ is Fredholm and use the decompositions in
Lemma~\ref{lem:descomposiciones-fredholm}. Define $S\colon Y\to X$ by
\[
S(r+c):=T_0^{-1}r,
\qquad
r\in\operatorname{Ran}T,\quad c\in C.
\]
The projections associated with the two direct sums and $T_0^{-1}$ are continuous, so
$S\in\mathcal L(Y,X)$. If $P_N$ and $P_C$ are the projections onto $\ker T$ and
$C$, respectively, then
\[
ST=I_X-P_N,
\qquad
TS=I_Y-P_C.
\]
Both projections have finite rank and are therefore compact. This proves
(ii).

Now assume (ii). For $x\in\ker T$, the first identity in
\eqref{eq:parametrix-atkinson} gives $x=K_Xx$. The unit ball of $\ker T$ is therefore the image under $K_X$ of a bounded set and is relatively compact.
The unit ball of an infinite-dimensional normed space cannot be
relatively compact; hence $\ker T$ is finite-dimensional.

Choose a closed complement $X_0$ of $\ker T$ as in the first part of
Lemma~\ref{lem:descomposiciones-fredholm}. We claim that there exists $c>0$ such that
\begin{equation}
\|Tx\|_Y\geq c\|x\|_X,
\qquad x\in X_0.
\label{eq:cota-inferior-complemento-fredholm}
\end{equation}
If this estimate were false, there would be a sequence $(x_\nu)\subseteq X_0$
with $\|x_\nu\|_X=1$ and $Tx_\nu\to0$. By compactness, a subsequence of
$(K_Xx_\nu)$ would converge. The identity
\[
x_\nu=STx_\nu+K_Xx_\nu
\]
would show that the same subsequence of $(x_\nu)$ converges to some $x\in X_0$.
Then $\|x\|_X=1$ and $Tx=0$, contradicting
$X_0\cap\ker T=\{0\}$. This proves
\eqref{eq:cota-inferior-complemento-fredholm}.

The estimate implies that the range is closed. Indeed, if
$Tx_\nu\to y$, write $x_\nu=n_\nu+x_\nu^0$ with
$n_\nu\in\ker T$ and $x_\nu^0\in X_0$. Then $Tx_\nu^0=Tx_\nu$ and
\[
\|x_\nu^0-x_\mu^0\|_X
\leq\frac{1}{c}\|Tx_\nu-Tx_\mu\|_Y.
\]
The sequence $(x_\nu^0)$ converges to an element $x^0\in X_0$, and the continuity
of $T$ gives $y=Tx^0$.

Let $Q\colon Y\to Y/\operatorname{Ran}T$ be the quotient map. Since the range is closed,
the codomain of $Q$ is Banach. Set $K_Y=I_Y-TS$. If $y=Tx$, then
\[
K_Yy=Tx-TSTx=T(I_X-ST)x=TK_Xx\in\operatorname{Ran}T.
\]
Consequently, $K_Y$ induces a compact operator
\[
\overline K_Y\colon Y/\operatorname{Ran}T
\longrightarrow Y/\operatorname{Ran}T,
\qquad
\overline K_Y[ y ]:=[K_Yy].
\]
To verify compactness, observe that every class of norm less than one
has a representative of norm less than two. The image of these representatives
under $K_Y$ is relatively compact in $Y$, and its projection under $Q$
remains relatively compact in the quotient.
The second identity in \eqref{eq:parametrix-atkinson} implies
\[
\overline K_Y[ y ]=[y-TSy]=[y],
\]
so the identity on the cokernel is compact. Its unit ball is
relatively compact, and the cokernel is finite-dimensional. Thus $T$ is
Fredholm.

When $X=Y$, the two identities in \eqref{eq:parametrix-atkinson} are equivalent to
$[S][T]=[I_X]=[T][S]$ in $\mathcal C(X)$, which is precisely the
invertibility of $[T]$.
\end{proof}

The preceding theorem appears in its Calkin algebra formulation in
\cite[Theorem~3.2, pp.~63--65]{BleeckerBoossIndex}. The proof reveals a
point that will be decisive in geometry: there is no need to construct an exact inverse;
a parametrix with compact remainders contains all the Fredholm information.

\begin{example}[Unilateral shifts]
Let $\ell^2(\mathbb N_0)$ be the space of square-summable sequences and define
\[
S_+(a_0,a_1,a_2,\ldots):=(0,a_0,a_1,\ldots),
\qquad
S_-(a_0,a_1,a_2,\ldots):=(a_1,a_2,a_3,\ldots).
\]
Both operators have norm one and
\[
S_-S_+=I,
\qquad
S_+S_-=I-\Pi_0,
\]
where $\Pi_0(a_0,a_1,\ldots)=(a_0,0,0,\ldots)$ has rank one. Thus each
is a parametrix for the other modulo compact operators. Moreover,
\[
\ker S_+=\{0\},
\qquad
\operatorname{Ran}S_+=\{(a_j)\in\ell^2(\mathbb N_0)\mid a_0=0\},
\]
whereas $\ker S_-=\operatorname{span}\{(1,0,0,\ldots)\}$ and $S_-$ is
surjective. Therefore,
\[
\operatorname{ind}S_+=-1,
\qquad
\operatorname{ind}S_-=1.
\]
This example shows that the index is not a formal difference between two
infinite dimensions: it records the finite defect left after inverting
the operator modulo compact operators.
\end{example}

\section{Stability and index calculus}

The definition of the index uses two dimensions that may change separately.
Their difference, however, remains fixed under perturbations that preserve
the Fredholm property.

\begin{theorem}[Stability of the Fredholm index]
\label{teo:estabilidad-indice-fredholm}
\index{Fredholm index!stability}
Let $X$ and $Y$ be Banach spaces.
\begin{enumerate}[label=(\alph*)]
\item The set
\[
\operatorname{Fred}(X,Y)
:=\{T\in\mathcal L(X,Y)\mid T\text{ is Fredholm}\}
\]
is open in the operator norm, and
$T\mapsto\operatorname{ind}T$ is locally constant.
\item If $T\in\operatorname{Fred}(X,Y)$ and $K\in\mathcal K(X,Y)$, then
$T+K$ is Fredholm and
\[
\operatorname{ind}(T+K)=\operatorname{ind}T.
\]
\item If $Z$ is a connected topological space and
$z\mapsto T_z\in\operatorname{Fred}(X,Y)$ is norm-continuous, then
$\operatorname{ind}T_z$ is independent of $z$.
\end{enumerate}
\end{theorem}

\begin{proof}
Fix $T\in\operatorname{Fred}(X,Y)$ and the decompositions
\[
X=X_0\oplus N,
\qquad
Y=R\oplus C,
\qquad
N=\ker T,\quad R=\operatorname{Ran}T,
\]
of Lemma~\ref{lem:descomposiciones-fredholm}. With respect to these, $T$ has
matrix
\[
T=
\begin{pmatrix}
T_0&0\\
0&0
\end{pmatrix},
\qquad
T_0\colon X_0\longrightarrow R
\text{ invertible}.
\]
For $A$ sufficiently close to $T$, write
\[
A=
\begin{pmatrix}
a&b\\
c&d
\end{pmatrix}.
\]
Let $\pi_R\colon Y\to R$ be the continuous projection associated
with the decomposition. Since $a=\pi_RA|_{X_0}$,
$\|a-T_0\|\leq\|\pi_R\|\|A-T\|$. Set
\[
 \varepsilon_T:=
 \frac{1}{2(1+\|T_0^{-1}\|)(1+\|\pi_R\|)}>0.
\]
If $\|A-T\|<\varepsilon_T$, then
\[
 \|T_0^{-1}(a-T_0)\|
 \leq\|T_0^{-1}\|\|\pi_R\|\|A-T\|
 <\frac{1}{2}.
\]
Thus $a=T_0(I_{X_0}+T_0^{-1}(a-T_0))$ is invertible by the Neumann
series. This criterion also applies when $X_0=R=\{0\}$, since the
unique operator from that space to itself is its identity and inverse.

Multiplying $A$ on the left and right by invertible operators yields
\[
\begin{pmatrix}
I_R&0\\
-ca^{-1}&I_C
\end{pmatrix}
\begin{pmatrix}
a&b\\
c&d
\end{pmatrix}
\begin{pmatrix}
I_{X_0}&-a^{-1}b\\
0&I_N
\end{pmatrix}
=
\begin{pmatrix}
a&0\\
0&d-ca^{-1}b
\end{pmatrix}.
\]
The second block acts between the finite-dimensional spaces $N$ and $C$.
Composition on the left or right with an isomorphism identifies the
corresponding kernels and cokernels and does not change the index. If
$D:=d-ca^{-1}b$, the rank--nullity theorem gives
\[
\dim\ker D-\dim(C/\operatorname{Ran}D)
=(\dim N-\dim\operatorname{Ran}D)
-(\dim C-\dim\operatorname{Ran}D)
=\dim N-\dim C.
\]
Thus $A$ is Fredholm and
\[
\operatorname{ind}A
=\dim N-\dim C
=\operatorname{ind}T.
\]
This proves both openness and local constancy.

If $K$ is compact and $S$ is a parametrix for $T$, then
\[
S(T+K)-I_X=(ST-I_X)+SK,
\qquad
(T+K)S-I_Y=(TS-I_Y)+KS
\]
are compact. Theorem~\ref{teo:atkinson-fredholm} shows that $T+tK$ is
Fredholm for every $t\in[0,1]$. Since $t\mapsto T+tK$ is continuous and $[0,1]$ is
connected, the local constancy just proved gives
$\operatorname{ind}(T+K)=\operatorname{ind}T$. The same argument, applied to
the continuous image of any connected space $Z$, proves (c).
\end{proof}

For local constancy, see also
\cite[Theorem~3.11, pp.~68--71]{BleeckerBoossIndex}. The block argument
explains why the dimensions of the kernel and cokernel need not
be constant: all their variation is confined to the finite-dimensional operator
$d-ca^{-1}b$, whose index is always $\dim N-\dim C$.

\begin{proposition}[Rules for index calculus]
\label{prop:calculo-indice-fredholm}
Let $T\colon X\to Y$ and $U\colon Y\to Z$ be Fredholm operators between Banach
spaces, and let $V\colon X_1\to Y_1$ be another Fredholm operator. Then
\begin{align}
\operatorname{ind}(T\oplus V)
&=\operatorname{ind}T+\operatorname{ind}V,
\label{eq:indice-suma-directa-fredholm}\\
\operatorname{ind}(UT)
&=\operatorname{ind}T+\operatorname{ind}U.
\label{eq:indice-composicion-fredholm}
\end{align}
If $X$ and $Y$ are Hilbert spaces and $T^*\colon Y\to X$ is the adjoint, then
$T^*$ is Fredholm and
\begin{equation}
\operatorname{ind}T^*=-\operatorname{ind}T.
\label{eq:indice-adjunto-fredholm}
\end{equation}
\end{proposition}

\begin{proof}
The direct sum satisfies
\[
\ker(T\oplus V)=\ker T\oplus\ker V,
\qquad
\operatorname{coker}(T\oplus V)
\cong\operatorname{coker}T\oplus\operatorname{coker}V,
\]
which proves \eqref{eq:indice-suma-directa-fredholm}.

Let $S_T\colon Y\to X$ and $S_U\colon Z\to Y$ be parametrices for $T$ and $U$.
Since compact operators form a two-sided ideal,
$S_TS_U$ is a parametrix for $UT$ modulo compact operators. By
Theorem~\ref{teo:atkinson-fredholm}, $UT$ is Fredholm. Consider the
sequence
\begin{equation}
\begin{aligned}
0&\longrightarrow\ker T
\longrightarrow\ker(UT)
\xrightarrow{\,T\,}\ker U
\xrightarrow{\,q_T\,}\operatorname{coker}T\\
&\xrightarrow{\,\overline U\,}\operatorname{coker}(UT)
\xrightarrow{\,q_U\,}\operatorname{coker}U
\longrightarrow0,
\end{aligned}
\label{eq:sucesion-exacta-composicion-fredholm}
\end{equation}
where $q_T(y)=[y]$, $\overline U([y])=[Uy]$, and $q_U([z])=[z]$. These maps
are well defined: $U(\operatorname{Ran}T)=\operatorname{Ran}(UT)$, and the
quotient by $\operatorname{Ran}(UT)$ projects naturally onto the quotient
by $\operatorname{Ran}U$. Let us verify exactness at the points where it is not
tautological. If $y\in\ker U$ and $q_T(y)=0$, then $y=Tx$ and
$x\in\ker(UT)$, so $y$ lies in the image of the third map. If
$[y]\in\operatorname{coker}T$ and $[Uy]=0$ in
$\operatorname{coker}(UT)$, there exists $x\in X$ such that $Uy=UTx$; therefore,
$y-Tx\in\ker U$ and $[y]=q_T(y-Tx)$. Finally, if
$[z]\in\operatorname{coker}(UT)$ vanishes in $\operatorname{coker}U$, then
$z=Uy$ and $[z]=\overline U([y])$. The sequence is exact.

All its terms are finite-dimensional. The alternating sum of the
dimensions in a finite exact sequence is zero; applying this to
\eqref{eq:sucesion-exacta-composicion-fredholm} gives
\[
\dim\ker T-\dim\ker(UT)+\dim\ker U
-\dim\operatorname{coker}T
+\dim\operatorname{coker}(UT)-\dim\operatorname{coker}U=0,
\]
which is equivalent to \eqref{eq:indice-composicion-fredholm}.

Now suppose that $X$ and $Y$ are Hilbert spaces. Their orthogonal decompositions
are
\[
X=\ker T\oplus(\ker T)^\perp,
\qquad
Y=\operatorname{Ran}T\oplus(\operatorname{Ran}T)^\perp.
\]
The restriction
$T_0\colon(\ker T)^\perp\to\operatorname{Ran}T$ is an isomorphism. Its adjoint
$T_0^*\colon\operatorname{Ran}T\to(\ker T)^\perp$ is also an isomorphism.
Moreover,
\[
\ker T^*=(\operatorname{Ran}T)^\perp,
\qquad
\operatorname{Ran}T^*=(\ker T)^\perp.
\]
The first equality follows directly from the definition of the adjoint;
the second follows from the invertibility of $T_0^*$ and the fact that $T^*$ vanishes on
$(\operatorname{Ran}T)^\perp$. In particular, $T^*$ has closed range,
\[
\dim\ker T^*=\dim\operatorname{coker}T,
\qquad
\dim\operatorname{coker}T^*=\dim\ker T,
\]
and we obtain \eqref{eq:indice-adjunto-fredholm}.
\end{proof}

\section{Continuous families and finite-dimensional stabilization}

For a norm-continuous family of Fredholm operators, the numerical index
is constant on each connected component of the parameter space. When the
parameter ranges over a compact space, a single finite-dimensional
correction can be chosen to make all the operators surjective simultaneously.
This will be the appropriate form for constructing a $K$--theory class in the next
chapter.

\begin{proposition}[Stabilization of a Fredholm family]
\label{prop:estabilizacion-familias-fredholm}
Let $B$ be a compact Hausdorff space and let
\[
T\colon B\longrightarrow\operatorname{Fred}(X,Y),
\qquad b\longmapsto T_b,
\]
be continuous in the operator norm. There exists a
finite-dimensional subspace $V\subseteq Y$ such that
\begin{equation}
\operatorname{Ran}T_b+V=Y,
\qquad b\in B.
\label{eq:transversal-finita-familia-fredholm}
\end{equation}
The surjective operators
\[
\widetilde T_b\colon X\oplus V\longrightarrow Y,
\qquad
\widetilde T_b(x,v):=T_bx+v,
\]
have kernels forming a vector bundle of locally constant rank over
$B$. If $V\subseteq V'$ are two stabilizations satisfying
\eqref{eq:transversal-finita-familia-fredholm}, then
\begin{equation}
[\ker\widetilde T^{\,V'}]-[B\times V']
=[\ker\widetilde T^{\,V}]-[B\times V]
\label{eq:independencia-estabilizacion-indice-familia}
\end{equation}
as stable differences of vector bundles.
\end{proposition}

\begin{proof}
Fix $b_0\in B$ and choose a finite-dimensional complement
$C_{b_0}$ of $\operatorname{Ran}T_{b_0}$. The operator
\[
A_{b_0}\colon X\oplus C_{b_0}\longrightarrow Y,
\qquad A_{b_0}(x,c)=T_{b_0}x+c,
\]
is surjective. Since
$C_{b_0}\cap\operatorname{Ran}T_{b_0}=\{0\}$, its kernel is
$\ker T_{b_0}\oplus\{0\}$ and is therefore finite-dimensional.
Hence the kernel of $A_{b_0}$ admits a closed complement $W_{b_0}$,
and $A_{b_0}\restriction_{W_{b_0}}\colon W_{b_0}\to Y$ is an isomorphism.
Write $R_{b_0}\colon Y\to W_{b_0}$ for its inverse and
$\pi_X\colon X\oplus C_{b_0}\to X$ for the projection. For $b\in B$, define
$A_b(x,c):=T_bx+c$ and
\[
 E_b:=R_{b_0}(T_b-T_{b_0})\pi_X\restriction_{W_{b_0}},
 \qquad
 \mathcal U_{b_0}:=\{b\in B\mid\|E_b\|_{\mathcal L(W_{b_0})}<\frac{1}{2}\}.
\]
Continuity of $T$ makes $\mathcal U_{b_0}$ an open subset of $B$ containing
$b_0$. For $b\in\mathcal U_{b_0}$,
\[
 A_b\restriction_{W_{b_0}}
 =(A_{b_0}\restriction_{W_{b_0}})(I_{W_{b_0}}+E_b),
 \qquad
 (A_b\restriction_{W_{b_0}})^{-1}
 =\sum_{j=0}^{\infty}(-E_b)^jR_{b_0}.
\]
The series converges in operator norm, and its norm is at most
$2\|R_{b_0}\|$. In particular, $A_b$ is surjective and
$\operatorname{Ran}T_b+C_{b_0}=Y$ throughout $\mathcal U_{b_0}$.

These neighborhoods cover $B$. Choose a finite subcover, associated with
$b_1,\ldots,b_N$, and set
\[
V:=C_{b_1}+\cdots+C_{b_N}.
\]
The space $V$ is finite-dimensional and satisfies
\eqref{eq:transversal-finita-familia-fredholm}.

To prove the bundle structure, fix $b_0$ again and set
$N_0:=\ker\widetilde T_{b_0}$. The projection $N_0\to V$ has kernel
$\ker T_{b_0}\oplus\{0\}$ and image contained in the finite-dimensional
space $V$; hence $N_0$ is finite-dimensional. Choose a closed complement
$W$ of $N_0$ in $X\oplus V$ and write
\[
 S_0:=(\widetilde T_{b_0}\restriction_W)^{-1}\colon Y\longrightarrow W
\]
for the inverse restriction. The set
\[
 \mathcal O_{b_0}:=
 \{b\in B\mid
   \|S_0(\widetilde T_b-\widetilde T_{b_0})\restriction_W\|_{\mathcal L(W)}<\frac{1}{2}\}
\]
is open and contains $b_0$. The same Neumann-series factorization gives,
for $b\in\mathcal O_{b_0}$, an inverse
\mbox{$Q_b:=(\widetilde T_b\restriction_W)^{-1}$}, continuous in operator norm.
For $n\in N_0$, define
\[
 \Phi_b(n):=n-Q_b\widetilde T_b n.
\]
Then $\Phi_b(n)\in\ker\widetilde T_b$. If
$z=n+w\in N_0\oplus W$ belongs to $\ker\widetilde T_b$, the equation
$\widetilde T_bz=0$ forces $w=-Q_b\widetilde T_bn$; thus $z=\Phi_b(n)$.
Let $\pi_{N_0}\colon X\oplus V\to N_0$ be the projection associated with
this decomposition. The maps
\[
 (b,n)\longmapsto(b,\Phi_b(n)),
 \qquad
 (b,z)\longmapsto(b,\pi_{N_0}z)
\]
are continuous inverses between $\mathcal O_{b_0}\times N_0$ and
$\{(b,z)\mid b\in\mathcal O_{b_0},\ \widetilde T_bz=0\}$, with the topology
induced by $B\times(X\oplus V)$. They are linear on each fiber.
This constructs the local trivializations and proves that the kernels
form a vector bundle of locally constant rank.

Finally, suppose that $V\subseteq V'$. The projection $V'\to V'/V$ induces an
exact sequence of bundles
\[
0\longrightarrow\ker\widetilde T^{\,V}
\longrightarrow\ker\widetilde T^{\,V'}
\longrightarrow B\times(V'/V)
\longrightarrow0.
\]
Surjectivity of the last arrow is checked fiber by fiber. Given
$w\in V'$, surjectivity of $\widetilde T_b^{\,V}$ provides
$(x,v)\in X\oplus V$ such that $T_bx+v=-w$; then $(x,v+w)$ belongs to
$\ker\widetilde T_b^{\,V'}$ and projects to $[w]\in V'/V$. A bundle
metric on the middle bundle, Euclidean if $\mathbb K=\mathbb R$ and
Hermitian if $\mathbb K=\mathbb C$, splits the sequence orthogonally, so
\[
\ker\widetilde T^{\,V'}
\cong\ker\widetilde T^{\,V}\oplus B\times(V'/V).
\]
Subtracting the trivial bundles $B\times V'$ and $B\times V$ yields
\eqref{eq:independencia-estabilizacion-indice-familia}.

It remains to compare two arbitrary stabilizations $V_0,V_1\subseteq Y$. The
subspace $W:=V_0+V_1$ is finite-dimensional and also satisfies
\eqref{eq:transversal-finita-familia-fredholm}. Applying the equality just
proved to the inclusions $V_0\subseteq W$ and $V_1\subseteq W$ yields
\[
 [\ker\widetilde T^{\,V_0}]-[B\times V_0]
 = [\ker\widetilde T^{\,W}]-[B\times W]
 = [\ker\widetilde T^{\,V_1}]-[B\times V_1].
\]
Thus the stable difference is independent of every choice
used in its construction.
\end{proof}

Stabilization replaces a family whose kernels may vary in dimension
with a well-defined difference of finite-dimensional bundles.
The Grothendieck group in the next chapter will allow us to consider this
difference without retaining the auxiliary subspace $V$. This is the
index bundle construction developed in
\cite[Theorem~3.30, pp.~82--85]{BleeckerBoossIndex}.

\section{Elliptic operators on closed manifolds}

Let $(M,\mathbf{g})$ be a closed smooth Riemannian manifold and let
$\mathbf{E},\mathbf{F}\longrightarrow M$ be finite-rank complex vector bundles equipped with Hermitian bundle metrics
and Hermitian connections. For $m,s\in\mathbb R$ and
\[
P\in\Psi^m_{\mathrm{cl}}(M;\mathbf{E},\mathbf{F}),
\],
Corollary~\ref{pseudo:cor-mapeo-sobolev-cerrada} gives a unique continuous
extension
\[
P_s\colon H^{s+m}(M,\mathbf{E})\longrightarrow H^s(M,\mathbf{F}).
\]
In what follows, $P$ will be elliptic.
Chapter~\ref{cap:operadores-pseudodiferenciales} constructed a parametrix
$Q\in\Psi^{-m}_{\mathrm{cl}}(M;\mathbf{F},\mathbf{E})$ with
\begin{equation}
QP=I_{\mathbf{E}}-R_{\mathbf{E}},\qquad PQ=I_{\mathbf{F}}-R_{\mathbf{F}},\qquad R_{\mathbf{E}},R_{\mathbf{F}}\in\Psi^{-\infty}.
\label{eq:parametrix-eliptica-indice-analitico}
\end{equation}

\begin{theorem}[Elliptic Fredholm operators and the analytic index]
\label{teo:fredholm-eliptico-indice-analitico}
\index{analytic index}
\index{elliptic operator!index}
Let $(M,\mathbf{g})$ be a closed Riemannian manifold, let
$\mathbf{E},\mathbf{F}\to M$ be finite-rank complex vector
bundles equipped with Hermitian bundle metrics, and let
$P\in\Psi^m_{\mathrm{cl}}(M;\mathbf{E},\mathbf{F})$ be elliptic.
For every $s\in\mathbb R$, the operator
\[
P_s\colon H^{s+m}(M,\mathbf{E})\longrightarrow H^s(M,\mathbf{F})
\]
is Fredholm. Its kernel consists of the smooth solutions of $P\mathbf{u}=0$, and its
cokernel is identified, by means of the orthogonal projection extended by duality,
with the space of smooth solutions
of $P^*\mathbf{v}=0$, where $P^*$ is the Hermitian formal adjoint. In particular,
\begin{equation}
\operatorname{ind}P_s
=\dim\ker\bigl(P\colon\Gamma(M,\mathbf{E})\to\Gamma(M,\mathbf{F})\bigr)
-\dim\ker\bigl(P^*\colon\Gamma(M,\mathbf{F})\to\Gamma(M,\mathbf{E})\bigr).
\label{eq:indice-analitico-nucleos-suaves}
\end{equation}
The right-hand side is independent of $s$.
\end{theorem}

\begin{proof}
The parametrix in \eqref{eq:parametrix-eliptica-indice-analitico} induces
\[
Q_s\colon H^s(M,\mathbf{F})\longrightarrow H^{s+m}(M,\mathbf{E}).
\]
By Proposition~\ref{pseudo:prop-regularizantes-compactos}, $R_{\mathbf{E}}$ and $R_{\mathbf{F}}$
act compactly on the Sobolev spaces occurring in the two
identities. Applying Theorem~\ref{teo:atkinson-fredholm} to $P_s$ and $Q_s$
proves that $P_s$ is Fredholm.

If $\mathbf{u}\in H^{s+m}(M,\mathbf{E})$ and $P_s\mathbf{u}=0$, the first parametrix identity gives
$\mathbf{u}=R_{\mathbf{E}}\mathbf{u}$. A smoothing operator takes every distribution to a smooth
section, so $\mathbf{u}\in\Gamma(M,\mathbf{E})$. Conversely, every smooth solution
belongs to $H^{s+m}(M,\mathbf{E})$ because $M$ is compact. Thus
\begin{equation}
\ker P_s=\{\mathbf{u}\in\Gamma(M,\mathbf{E})\mid P\mathbf{u}=0\},
\label{eq:nucleo-eliptico-independiente-s}
\end{equation}
independently of $s$.

The formal adjoint $P^*\in\Psi^m_{\mathrm{cl}}(M;\mathbf{F},\mathbf{E})$ is characterized by
\begin{equation}
\int_M \mathbf{h}_{\mathbf{F}}(P\mathbf{u},\mathbf{v})\,d\lambda_{\mathbf{g}}
=\int_M \mathbf{h}_{\mathbf{E}}(\mathbf{u},P^*\mathbf{v})\,d\lambda_{\mathbf{g}},
\qquad
\mathbf{u}\in\Gamma(M,\mathbf{E}),\quad \mathbf{v}\in\Gamma(M,\mathbf{F}).
\label{eq:adjunto-formal-indice-analitico}
\end{equation}
The principal symbol of $P^*$ is the adjoint of the principal symbol of $P$ by
Theorem~\ref{teo:simbolo-adjunto-pseudo}; therefore, $P^*$ is elliptic.
Proposition~\ref{prop:dualidad-sobolev-haces-geometria-acotada}, applied to
the closed manifold, identifies the continuous annihilator of
$\operatorname{Ran}P_s\subseteq H^s(M,\mathbf{F})$ with
\[
\ker\left(
P^*\colon H^{-s}(M,\mathbf{F})\longrightarrow H^{-s-m}(M,\mathbf{E})
\right),
\]
by an antilinear identification, since the Hermitian inner product is
linear in its first variable. The elliptic regularity of
Corollary~\ref{pseudo:cor-regularidad-eliptica-sobolev-local} shows that this
kernel consists of smooth sections. Denote it by $K$ and choose an
orthonormal basis $\mathbf{v}_1,\ldots,\mathbf{v}_q$ of $K$ for the
$L^2(M,\mathbf{F})$ inner product. For $\mathbf{f}\in H^s(M,\mathbf{F})$ define
\[
 \Pi_K\mathbf{f}
 :=\sum_{j=1}^q\langle\mathbf{f},\mathbf{v}_j\rangle\mathbf{v}_j.
\]
Each pairing is continuous on $H^s(M,\mathbf F)$, because $\mathbf{v}_j$ is smooth
and belongs to $H^{-s}(M,\mathbf F)$. Thus $\Pi_K:H^s(M,\mathbf{F})\to K$ is
linear and continuous, agrees with the orthogonal projection when
$\mathbf{f}\in L^2(M,\mathbf F)\cap H^s(M,\mathbf F)$, and is the identity on $K$.
Its kernel is $\operatorname{Ran}P_s$: inclusion of the range follows
from the formal adjoint; conversely, if all the pairings vanish,
all continuous functionals annihilating $\operatorname{Ran}P_s$
vanish at $\mathbf{f}$. Hahn--Banach and closedness of the range imply
$\mathbf{f}\in\operatorname{Ran}P_s$. Thus $\Pi_K$ induces a complex-linear
isomorphism $\operatorname{coker}P_s\cong K$, independent of the chosen
orthonormal basis.

Since $\operatorname{Ran}P_s$ is closed, the dual of the cokernel is this annihilator.
Both are finite-dimensional and therefore have the same dimension. Together
with \eqref{eq:nucleo-eliptico-independiente-s}, this gives
\eqref{eq:indice-analitico-nucleos-suaves} and proves independence of $s$.
\end{proof}

\begin{definition}[Analytic index]
\label{def:indice-analitico}
Let $(M,\mathbf{g})$ be a closed Riemannian manifold, let
$\mathbf{E},\mathbf{F}\to M$ be finite-rank complex vector
bundles equipped with Hermitian bundle metrics, and let
$P\in\Psi^m_{\mathrm{cl}}(M;\mathbf{E},\mathbf{F})$ be elliptic.
The \textbf{analytic index} of a classical elliptic operator
$P$ is the integer
\[
\operatorname{ind}_{\mathrm a}(P)
:=\operatorname{ind}\left(
P_s\colon H^{s+m}(M,\mathbf{E})\longrightarrow H^s(M,\mathbf{F})
\right),
\]
where $s\in\mathbb R$ is arbitrary.
Theorem~\ref{teo:fredholm-eliptico-indice-analitico} ensures that the definition
is independent of $s$.
\end{definition}

Formula \eqref{eq:indice-analitico-nucleos-suaves} can be written as
\[
\operatorname{ind}_{\mathrm a}(P)=\dim\ker P-\dim\ker P^*.
\]
Here both kernels are taken among smooth sections. Using the Hilbert space adjoint
of $P_s$ for specific Sobolev inner products gives an operator different from
$P^*$; the Riesz isomorphism identifies their kernels and yields the same
dimension. The formal adjoint preserves the geometric content and the
principal symbol.

\begin{corollary}[Properties of the analytic index]
\label{cor:propiedades-indice-analitico}
Let $M$ be a closed smooth manifold, let
$\mathbf{E},\mathbf{F},\mathbf{E}_0,\mathbf{F}_0,\mathbf{G}\to M$ be finite-rank
smooth complex vector bundles, and let
$P\in\Psi^m_{\mathrm{cl}}(M;\mathbf{E},\mathbf{F})$,
$P_0\in\Psi^m_{\mathrm{cl}}(M;\mathbf{E}_0,\mathbf{F}_0)$, and
$Q\in\Psi^\ell_{\mathrm{cl}}(M;\mathbf{F},\mathbf{G})$ be elliptic. Then
\begin{align*}
\operatorname{ind}_{\mathrm a}(P\oplus P_0)
&=\operatorname{ind}_{\mathrm a}(P)+\operatorname{ind}_{\mathrm a}(P_0),\\
\operatorname{ind}_{\mathrm a}(QP)
&=\operatorname{ind}_{\mathrm a}(P)+\operatorname{ind}_{\mathrm a}(Q),\\
\operatorname{ind}_{\mathrm a}(P^*)
&=-\operatorname{ind}_{\mathrm a}(P).
\end{align*}
If $P$ is invertible, or if $P=P^*$, then
$\operatorname{ind}_{\mathrm a}(P)=0$.
\end{corollary}

\begin{proof}
Choose Sobolev exponents making all compositions continuous. The
first two identities are, respectively,
\eqref{eq:indice-suma-directa-fredholm} and
\eqref{eq:indice-composicion-fredholm}; independence of the exponent allows us to
change scales without changing the integer. An invertible operator has zero kernel and
cokernel. For the third identity, apply
\eqref{eq:indice-analitico-nucleos-suaves} to $P$ and $P^*$:
since $(P^*)^*=P$, the two dimensions exchange places.
If $P=P^*$, this identity gives
$\operatorname{ind}_{\mathrm a}(P)=-\operatorname{ind}_{\mathrm a}(P)$, and since
the index is an integer, it must be zero.
\end{proof}

\section{Dependence on the principal symbol}

Let $\pi\colon T^*M\to M$ be the projection. The principal symbol of $P$ is a
homogeneous homomorphism
\[
\boldsymbol{\sigma}_m^\Psi(P)\colon\pi^*\mathbf{E}\longrightarrow\pi^*\mathbf{F}
\]
that is invertible over $T^*M\setminus0_M$. Its restriction to the unit cosphere bundle
$S^*M$ contains the same information, since homogeneity recovers the
symbol on every nonzero ray.

\begin{proposition}[Continuous quantization of a homotopy of symbols]
\label{prop:cuantizacion-homotopia-simbolos}
Let $M$ be a closed smooth manifold, let
$\mathbf{E},\mathbf{F}\to M$ be finite-rank smooth complex vector
bundles, and let $a_t\colon\pi^*\mathbf{E}\to\pi^*\mathbf{F}$, $t\in[0,1]$, be classical principal symbols of
order $m$, depending continuously on $t$ in the $C^\infty$ topology on
$S^*M$. There exists a family
\[
P_t\in\Psi^m_{\mathrm{cl}}(M;\mathbf{E},\mathbf{F})
\]
continuous as a map
\[
[0,1]\longrightarrow
\mathcal L\bigl(H^{s+m}(M,\mathbf{E}),H^s(M,\mathbf{F})\bigr)
\]
for every $s\in\mathbb R$, such that $\boldsymbol{\sigma}_m^\Psi(P_t)=a_t$. If each $a_t$ is
invertible, each $P_t$ is elliptic.
\end{proposition}

\begin{proof}
Choose a finite atlas, bundle frames, a partition of unity, and cutoffs
in the position and frequency variables as in the global construction of
Chapter~\ref{cap:operadores-pseudodiferenciales}. Extend each restriction
homogeneously from $S^*M$ for $\|\xi\|_{\mathbf{g}}\geq1$ and multiply it by a
cutoff vanishing near the zero section. The localized sum of the
quantizations has principal symbol $a_t$; terms arising from
changes of chart and cutoffs have order at most $m-1$ and do not alter this
symbol.

In a finite atlas, every symbol seminorm of the full symbol of $P_t-P_{t_0}$
is bounded by a finite sum of $C^N$ seminorms of $a_t-a_{t_0}$ on
$S^*M$. Theorem~\ref{pseudo:teo-mapeo-sobolev-euclidiano}, followed by the
finite localization used in
Corollary~\ref{pseudo:cor-mapeo-sobolev-cerrada}, shows that
\[
\|P_t-P_{t_0}\|_{\mathcal L(H^{s+m}(M,\mathbf{E}),H^s(M,\mathbf{F}))}
\longrightarrow0
\qquad(t\longrightarrow t_0).
\]
The last assertion is the definition of ellipticity.
\end{proof}

\begin{corollary}[The index depends only on the symbol]
\label{cor:indice-analitico-simbolo}
Let $M$ be a closed smooth manifold and let
$\mathbf{E},\mathbf{F}\to M$ be finite-rank smooth complex vector
bundles.
Let $P_0,P_1\in\Psi^m_{\mathrm{cl}}(M;\mathbf{E},\mathbf{F})$ be elliptic.
\begin{enumerate}[label=(\alph*)]
\item If $\boldsymbol{\sigma}_m^\Psi(P_0)=\boldsymbol{\sigma}_m^\Psi(P_1)$, then
\[
\operatorname{ind}_{\mathrm a}(P_0)
=\operatorname{ind}_{\mathrm a}(P_1).
\]
\item If the restrictions of $\boldsymbol{\sigma}_m^\Psi(P_0)$ and $\boldsymbol{\sigma}_m^\Psi(P_1)$ to $S^*M$ are
joined by a homotopy of bundle isomorphisms continuous in the
$C^\infty$ topology on sections, then the same equality holds.
\end{enumerate}
\end{corollary}

\begin{proof}
In (a), the exact sequence in
Theorem~\ref{teo:sucesion-exacta-simbolo-principal-pseudo} shows that the
difference $P_1-P_0$ has order at most $m-1$. For any
$s\in\mathbb R$, it factors as
\[
H^{s+m}(M,\mathbf{E})
\xrightarrow{\,P_1-P_0\,}
H^{s+1}(M,\mathbf{F})
\hookrightarrow H^s(M,\mathbf{F}).
\]
The first arrow is continuous by the pseudodifferential mapping theorem and the
second is compact by Rellich--Kondrashov. Thus $P_1-P_0$ is compact
between the spaces defining the index, and
Theorem~\ref{teo:estabilidad-indice-fredholm} proves the equality.

In (b), Proposition~\ref{prop:cuantizacion-homotopia-simbolos} produces a
norm-continuous family $(\widehat P_t)_{t\in[0,1]}$ of elliptic operators with
the prescribed symbols. It is a Fredholm family by
Theorem~\ref{teo:fredholm-eliptico-indice-analitico}, so part (c)
of Theorem~\ref{teo:estabilidad-indice-fredholm} gives
\[
\operatorname{ind}_{\mathrm a}(\widehat P_0)
=\operatorname{ind}_{\mathrm a}(\widehat P_1).
\]
Since $\widehat P_j$ and $P_j$ have the same principal symbol for $j=0,1$,
part (a) completes the proof.
\end{proof}

To see directly that small homotopies remain elliptic, fix
metrics $\mathbf{h}_{\mathbf{E}}$ and $\mathbf{h}_{\mathbf{F}}$ on the bundles. Compactness of $S^*M$ gives
\[
\delta
:=\inf_{(x,\xi)\in S^*M}
\left\|\boldsymbol{\sigma}_m^\Psi(P)(x,\xi)^{-1}\right\|_{
\mathcal L((\mathbf{F}_x,\mathbf{h}_{\mathbf{F}}),(\mathbf{E}_x,\mathbf{h}_{\mathbf{E}}))}^{-1}>0.
\]
If another symbol $b_m$ satisfies
\[
\sup_{(x,\xi)\in S^*M}
\|b_m(x,\xi)-\boldsymbol{\sigma}_m^\Psi(P)(x,\xi)\|_{
\mathcal L((\mathbf{E}_x,\mathbf{h}_{\mathbf{E}}),(\mathbf{F}_x,\mathbf{h}_{\mathbf{F}}))}<\delta,
\]
then
\[
b_m(x,\xi)
=\boldsymbol{\sigma}_m^\Psi(P)(x,\xi)
\left[I+
\boldsymbol{\sigma}_m^\Psi(P)(x,\xi)^{-1}
\bigl(b_m(x,\xi)-\boldsymbol{\sigma}_m^\Psi(P)(x,\xi)\bigr)
\right]
\]
is invertible at every point, because the second factor can be inverted using the
Neumann series. The same calculation applies to the segment joining the two symbols.

Corollary~\ref{cor:indice-analitico-simbolo} thus reduces the analytic problem to
a topological one: turning the isomorphism
$\boldsymbol{\sigma}_m^\Psi(P)$ over $T^*M\setminus0_M$ into a stable class and proving that the
index defines a homomorphism on these classes. This will be the construction in
Chapter~\ref{cap:k-teoria-clase-simbolo}.

\chapter{K-theory of vector bundles and the symbol class}
\label{cap:k-teoria-clase-simbolo}

Let $X$ be a compact Hausdorff space. Direct sum of complex vector
bundles over $X$ gives a commutative monoid, but does not
allow bundles to be subtracted. The Grothendieck construction remedies precisely
this deficiency: it adds formal inverses without identifying any more elements
than direct sum requires. In this sense the construction is
algebraic. The resulting group, however, depends on the base space and on
homotopies of maps; it is therefore called \emph{topological
$K$--theory}. It is distinct from the algebraic $K$--theory of rings.

The usefulness of this construction for an elliptic operator is already apparent
in its principal symbol. If $M$ is closed and
\[
 P\in\Psi^m_{\mathrm{cl}}(M;\mathbf{E},\mathbf{F})
\]
is elliptic, then
\[
 \boldsymbol{\sigma}_m^\Psi(P)\colon \pi^*\mathbf{E}\longrightarrow\pi^*\mathbf{F}
\]
is an isomorphism away from the zero section of $T^*M$. The pair of bundles and this
isomorphism determine a compactly supported class on $T^*M$. In this
chapter we construct this class and prove the properties of the operations
used to define the topological and analytic indices. The construction follows
\cite[Chapter~10, in particular pp.~262--269]{BleeckerBoossIndex}.

All bundles considered will be complex, finite-rank, and locally
trivial. Over a smooth manifold they will also be smooth. Real $K$--theory
requires different groups and will not enter the proof of
the complex index theorem.

\begin{semblanzaHistorica}{Grothendieck and Atiyah: subtracting bundles}
Grothendieck completion turns a monoid of geometric objects into a group of formal differences. Applied to vector bundles, this idea gives a stable expression for the difference between the domain and codomain of an elliptic symbol. Atiyah developed topological $K$--theory as the natural language for these classes and, together with Singer, made it the topological side of the index theorem. Stabilization does not erase the geometry: it retains precisely the information that persists when trivial bundles are added.
\end{semblanzaHistorica}

\section{Grothendieck completion}
\label{sec:grothendieck-haces}

Denote by $\operatorname{Vect}_{\mathbb C}(X)$ the set of isomorphism
classes of complex vector bundles over $X$. The operation
\[
 [\mathbf{E}]+[\mathbf{F}]:=[\mathbf{E}\oplus \mathbf{F}]
\]
makes it a commutative monoid whose identity element is the class of the
rank-zero bundle. Cancellation does not hold in general:
$\mathbf{E}\oplus \mathbf{G}\cong \mathbf{F}\oplus \mathbf{G}$ does not necessarily imply $\mathbf{E}\cong \mathbf{F}$.
This is why differences of bundles must not be defined
by first imposing a cancellation property that does not hold.

\begin{theorem}[Grothendieck group]
\label{teo:grupo-grothendieck-monoide}
Let $(A,+,0)$ be a commutative monoid. On $A\times A$ define
\begin{equation}
 (a,b)\sim(c,d)
 \quad\Longleftrightarrow\quad
 \text{there exists }e\in A\text{ such that }a+d+e=c+b+e.
 \label{eq:equivalencia-grothendieck}
\end{equation}
Then
\[
 G(A):=(A\times A)/{\sim}
\]
is an abelian group with
\[
 [(a,b)]+[(c,d)]=[(a+c,b+d)],
 \qquad
 -[(a,b)]=[(b,a)].
\]
The map
\[
 \iota_A\colon A\longrightarrow G(A),
 \qquad
 \iota_A(a)=[(a,0)],
\]
is a monoid homomorphism and satisfies the following universal
property: if $H$ is an abelian group and $h\colon A\to H$ is a
monoid homomorphism, there exists a unique group homomorphism
$\widetilde h\colon G(A)\to H$ such that
$\widetilde h\circ\iota_A=h$.
\end{theorem}

\begin{proof}
Let us prove that this is an equivalence relation. Reflexivity follows by taking
$e=0$, and symmetry by interchanging the two sides. If
\[
 a+d+e=c+b+e
 \quad\text{and}\quad
 c+f+e'=g+d+e',
\]
then adding the equalities and regrouping yields
\[
 a+f+(c+d+e+e')=g+b+(c+d+e+e'),
\]
so $(a,b)\sim(g,f)$. Thus the relation is transitive.

If $(a,b)$ and $(c,d)$ are replaced by equivalent representatives,
adding the witnesses for the two equivalences proves that
$(a+c,b+d)$ remains in the same class. Thus addition in $G(A)$ is well defined.
Moreover,
\[
 [(a,b)]+[(b,a)]=[(a+b,a+b)]=[(0,0)],
\]
so each element has the stated inverse.

For the universal property, define
\[
 \widetilde h([(a,b)]):=h(a)-h(b).
\]
If $(a,b)\sim(c,d)$, there exists $e\in A$ with
$a+d+e=c+b+e$; applying $h$ and canceling in the group $H$ gives
$h(a)-h(b)=h(c)-h(d)$. Consequently, $\widetilde h$ is well
defined. It is a homomorphism, factors $h$, and is unique because every class
can be written as
\[
 [(a,b)]=\iota_A(a)-\iota_A(b).
\]
\end{proof}

\begin{remark}
The relation is often defined by
\[
 (a,b)\sim(c,d)
 \quad\Longleftrightarrow\quad
 \text{there exist }u,v\in A\text{ with }
 (a+u,b+u)=(c+v,d+v).
\]
This form and \eqref{eq:equivalencia-grothendieck} are equivalent: in one
direction take $e=u+v$, and in the other take $u=d+e$ and $v=b+e$.
The second form is the one used in
\cite[Theorem~10.5, pp.~262--263]{BleeckerBoossIndex}.
\end{remark}

\begin{definition}[$K^0$ of a compact space]
\label{def:k0-espacio-compacto}
\label{def:k-teoria-espacio-compacto}
Let $X$ be compact and Hausdorff. Its complex $K$--theory group is
\[
 K^0(X):=G\bigl(\operatorname{Vect}_{\mathbb C}(X)\bigr).
\]
The class of the pair $(\mathbf{E},\mathbf{F})$ is denoted by
\[
 [\mathbf{E}]-[\mathbf{F}]\in K^0(X).
\]
If $\underline{\mathbb C}^{\,r}_X=X\times\mathbb C^r$ is the trivial bundle,
we abbreviate its class by $r$ when the base space is unambiguous.
\end{definition}

Equality between differences has a stable description that we will use
repeatedly.

\begin{proposition}[Criterion for stable equality]
\label{prop:igualdad-estable-k0}
Let $\mathbf{E},\mathbf{F},\mathbf{E}',\mathbf{F}'$ be complex vector bundles over a compact
Hausdorff space $X$. Then
\begin{equation}
 [\mathbf{E}]-[\mathbf{F}]=[\mathbf{E}']-[\mathbf{F}']
 \label{eq:igualdad-diferencias-k}
\end{equation}
if and only if there exists a bundle $\mathbf{G}\to X$ such that
\begin{equation}
 \mathbf{E}\oplus \mathbf{F}'\oplus \mathbf{G}\cong \mathbf{E}'\oplus \mathbf{F}\oplus \mathbf{G}.
 \label{eq:isomorfismo-estable-haz-g}
\end{equation}
This is also equivalent to the existence of $N\in\mathbb N_0$ for which
\begin{equation}
 \mathbf{E}\oplus \mathbf{F}'\oplus\underline{\mathbb C}^{\,N}_X
 \cong
 \mathbf{E}'\oplus \mathbf{F}\oplus\underline{\mathbb C}^{\,N}_X.
 \label{eq:isomorfismo-estable-trivial}
\end{equation}
Moreover, every class in $K^0(X)$ can be written as
\begin{equation}
 [\mathbf{V}]-[\underline{\mathbb C}^{\,N}_X]
 \label{eq:clase-k-como-haz-menos-trivial}
\end{equation}
for some bundle $\mathbf{V}\to X$ and some $N\in\mathbb N_0$.
\end{proposition}

\begin{proof}
The equivalence between \eqref{eq:igualdad-diferencias-k} and
\eqref{eq:isomorfismo-estable-haz-g} is precisely
\eqref{eq:equivalencia-grothendieck} applied to the monoid of bundle classes.

By the complement theorem for bundles over compact bases, stated
in Appendix~\ref{ap:topologia}, there exist a bundle $\mathbf{G}^\perp$ and an integer $N$
such that
\[
 \mathbf{G}\oplus \mathbf{G}^\perp\cong\underline{\mathbb C}^{\,N}_X.
\]
Adding $\mathbf{G}^\perp$ to both sides of
\eqref{eq:isomorfismo-estable-haz-g} yields
\eqref{eq:isomorfismo-estable-trivial}. The reverse implication follows
immediately by taking $\mathbf{G}$ trivial.

Finally, given $[\mathbf{E}]-[\mathbf{F}]$, choose a complement $\mathbf{F}^\perp$ with
$\mathbf{F}\oplus \mathbf{F}^\perp\cong\underline{\mathbb C}^{\,N}_X$. Then
\[
 [\mathbf{E}]-[\mathbf{F}]
 =[\mathbf{E}\oplus \mathbf{F}^\perp]-[\underline{\mathbb C}^{\,N}_X],
\]
which has the form \eqref{eq:clase-k-como-haz-menos-trivial}.
\end{proof}

The complement theorem is not a formal property of monoids. Its
proof uses a finite trivializing cover, a partition of
unity, and an inclusion of the bundle into a trivial bundle; the necessary
topological statements are collected in Appendix~\ref{ap:topologia}. The algebraic
consequence follows by applying that theorem in the preceding proposition.

\subsection{Product, pullback, and homotopy}

Tensor product distributes over direct sum. By the
universal property of the Grothendieck group, it induces a bilinear product
\begin{equation}
 \begin{aligned}
 K^0(X)\times K^0(X)&\longrightarrow K^0(X),\\
 \bigl([\mathbf{E}]-[\mathbf{F}],[\mathbf{G}]-[\mathbf{H}]\bigr)
 &\longmapsto
 [\mathbf{E}\otimes \mathbf{G}]+[\mathbf{F}\otimes \mathbf{H}]
 -[\mathbf{E}\otimes \mathbf{H}]-[\mathbf{F}\otimes \mathbf{G}].
 \end{aligned}
 \label{eq:producto-interno-k0}
\end{equation}
Associativity, commutativity, and distributivity follow from the
corresponding canonical bundle isomorphisms. The class of the
trivial rank-one bundle is the identity. Consequently, $K^0(X)$ is a commutative
ring.

If $f\colon Y\to X$ is continuous between compact Hausdorff spaces,
pullback satisfies
\[
 f^*(\mathbf{E}\oplus \mathbf{F})\cong f^*\mathbf{E}\oplus f^*\mathbf{F},
 \qquad
 f^*(\mathbf{E}\otimes \mathbf{F})\cong f^*\mathbf{E}\otimes f^*\mathbf{F}.
\]
It therefore defines a ring homomorphism
\[
 f^*\colon K^0(X)\longrightarrow K^0(Y),
 \qquad
 f^*([\mathbf{E}]-[\mathbf{F}])=[f^*\mathbf{E}]-[f^*\mathbf{F}].
\]
The identities
\[
 (f\circ g)^*=g^*\circ f^*,
 \qquad
 \operatorname{Id}_X^*=\operatorname{Id}_{K^0(X)}
\]
are verified first for bundles and then for differences. Thus $K^0$ is a
contravariant functor.

\begin{proposition}[Homotopy invariance]
\label{prop:homotopia-k0}
If $f_0,f_1\colon Y\to X$ are homotopic maps between compact Hausdorff
spaces, then
\[
 f_0^*=f_1^*\colon K^0(X)\longrightarrow K^0(Y).
\]
In particular, a homotopy equivalence induces an isomorphism in
$K^0$.
\end{proposition}

\begin{proof}
Let $H\colon Y\times[0,1]\to X$ be a homotopy from $f_0$ to $f_1$. The homotopy
invariance theorem for vector bundles in
Appendix~\ref{ap:topologia} gives
\[
 (H^*\mathbf{E})|_{Y\times\{0\}}\cong(H^*\mathbf{E})|_{Y\times\{1\}}
\]
for every bundle $\mathbf{E}\to X$. Consequently,
$f_0^*[\mathbf{E}]=f_1^*[\mathbf{E}]$. Since bundle classes generate $K^0(X)$, the two
homomorphisms agree on the whole group.
\end{proof}

\section{Reduced groups, pairs, and compact support}
\label{sec:k-pares-soporte-compacto}

Let $(X,x_0)$ be a pointed compact Hausdorff space. The inclusion
$i\colon\{x_0\}\hookrightarrow X$ induces the rank homomorphism at the
base point
\[
 i^*\colon K^0(X)\longrightarrow K^0(\{x_0\})\cong\mathbb Z.
\]

\begin{definition}[Reduced $K$--theory]
\label{def:k-teoria-reducida}
Define
\[
 \widetilde K^0(X):=\ker\bigl(i^*\colon K^0(X)\to\mathbb Z\bigr).
\]
Since the constant map $r\colon X\to\{x_0\}$ satisfies $i^*r^*=I$, we
obtain the natural decomposition
\[
 K^0(X)\cong\widetilde K^0(X)\oplus\mathbb Z.
\]
\end{definition}

Now let $A\subseteq X$ be a nonempty closed subspace and let $X/A$ be the
pointed space obtained by collapsing $A$ to a point.

\begin{definition}[Relative $K$--theory]
\label{def:k-teoria-relativa}
The relative group of the pair $(X,A)$ is
\[
 K^0(X,A):=\widetilde K^0(X/A).
\]
If $A=\varnothing$, we adopt the convention
$K^0(X,\varnothing)=K^0(X)$.
\end{definition}

The description by triples exhibits the symbol of an operator and is
valid for a closed subset of a compact Hausdorff space. In particular,
it covers pairs consisting of a disk bundle and its sphere bundle, as
well as their products with compact Hausdorff parameter spaces.

\begin{definition}[Relative triple]
\label{def:triple-relativo-k}
A triple over $(X,A)$ consists of two complex bundles $\mathbf{E},\mathbf{F}\to X$ and an
isomorphism
\[
 \alpha\colon \mathbf{E}|_A\longrightarrow \mathbf{F}|_A.
\]
We write $(\mathbf{E},\mathbf{F},\alpha)$. Direct sum is defined by
\[
 (\mathbf{E},\mathbf{F},\alpha)\oplus(\mathbf{E}',\mathbf{F}',\alpha')
 :=(\mathbf{E}\oplus \mathbf{E}',\mathbf{F}\oplus \mathbf{F}',\alpha\oplus\alpha').
\]
The group of triples is obtained by imposing the following relations:
\begin{enumerate}[label=(\roman*)]
 \item isomorphic triples represent the same class;
 \item a homotopy from $\alpha$ to $\alpha'$ through isomorphisms over
 $A$ does not change the class;
 \item if $\alpha$ extends to an isomorphism $\mathbf{E}\to \mathbf{F}$ over $X$, then
 $(\mathbf{E},\mathbf{F},\alpha)$ represents zero;
 \item if
 $\alpha\colon \mathbf{E}|_A\to \mathbf{F}|_A$ and
 $\beta\colon \mathbf{F}|_A\to \mathbf{G}|_A$, then
 \begin{equation}
 [(\mathbf{E},\mathbf{F},\alpha)]+[(\mathbf{F},\mathbf{G},\beta)]
 =[(\mathbf{E},\mathbf{G},\beta\circ\alpha)].
 \label{eq:relacion-composicion-triples}
 \end{equation}
\end{enumerate}
\end{definition}

The last relation contains both the existence of inverses and the
correct cancellation of an intermediate bundle. Indeed,
\[
 -[(\mathbf{E},\mathbf{F},\alpha)]=[(\mathbf{F},\mathbf{E},\alpha^{-1})]
\]
because the sum of the two classes is
$[(\mathbf{E},\mathbf{E},I_{\mathbf{E}|_A})]=0$. Addition is commutative: the isomorphism interchanging the
summands extends over all of $X$, so the interchange does not alter the
class. Thus the preceding relations produce an abelian group.

\begin{theorem}[The model by triples]
\label{teo:modelo-triples-k-relativa}
Let $X$ be a compact Hausdorff space and let $A\subseteq X$ be closed.
The group of triples in Definition~\ref{def:triple-relativo-k} is
naturally isomorphic to $K^0(X,A)$. If $A\neq\varnothing$, addition of
triples corresponds, under this isomorphism, to addition in
$\widetilde K^0(X/A)$.
\end{theorem}

\begin{proof}
If $A=\varnothing$, a triple contains only the two bundles, and the relations
are those of Grothendieck completion. Thus suppose that
$A\neq\varnothing$. We will use stable complements, gluing,
homotopy invariance, and exactness of $K$--theory from
Appendix~\ref{ap:topologia}, all in the category of compact
Hausdorff spaces. In the degree-one terms below we use
$K^1(Z)=\widetilde K^0(\Sigma Z_+)$, where $Z_+$ adjoins a disjoint base
point and $\Sigma$ is reduced suspension; this definition is recorded
explicitly later in this section.

First recall why gluing along a closed subset is possible without assuming
that the pair is a CW pair. An isomorphism
$\alpha:\mathbf E|_A\to\mathbf F|_A$ extends to an isomorphism on
some neighborhood of $A$. Indeed, the bundle
$\operatorname{Hom}(\mathbf E,\mathbf F)$ is a summand of a finite-rank
trivial bundle by Theorem~\ref{teo:complemento-haz-complejo-apendice}.
Regard $\alpha$ as a section of this trivial bundle over $A$, extend
its components to $X$ by the Tietze theorem, and project onto the
summand $\operatorname{Hom}(\mathbf E,\mathbf F)$. This gives a
homomorphism agreeing with $\alpha$ on $A$. Equality of the ranks
holds on a neighborhood of $A$, and invertibility is open among
matrices of that rank, so the homomorphism is invertible after
shrinking the neighborhood. Applied to a trivialization, this argument
also allows a bundle to descend to the quotient; see
\cite[Proposition~2.9, pp.~51--52]{HatcherVBKT}.

Form the double
\[
 D:=X_+\cup_A X_-,
\]
where $X_+$ and $X_-$ are two copies of $X$ and their copies of
$A$ are identified pointwise. It is compact and Hausdorff: the equivalence
relation is closed in the compact Hausdorff space $X_+\coprod X_-$. Denote by
$r:D\to X$ the map that is the identity on both copies, and by
$i_-:X\to D$ the inclusion of the second copy. We have
$r\circ i_-=\operatorname{Id}_X$ and a natural homeomorphism
$D/X_-\cong X/A$. Let $q:D\to X/A$ be the corresponding quotient map.

The exact sequence in
Theorem~\ref{teo:exactitud-excision-k-apendice}, applied to $(D,X_-)$,
contains
\[
 K^1(D)\xrightarrow{i_-^*}K^1(X)
 \xrightarrow{\partial}\widetilde K^0(X/A)
 \xrightarrow{q^*}K^0(D)\xrightarrow{i_-^*}K^0(X).
\]
The first map is surjective because $i_-^*r^*=I$ also in degree
one. By exactness, $\partial=0$, the map $q^*$ is injective, and
\[
 q^*:\widetilde K^0(X/A)\xrightarrow{\;\cong\;}
 \ker\bigl(i_-^*:K^0(D)\to K^0(X)\bigr).
\]
Thus it suffices to identify the group of triples with this last kernel.

Given $(\mathbf E,\mathbf F,\alpha)$, glue $\mathbf E$ over $X_+$
to $\mathbf F$ over $X_-$ by means of $\alpha$, and call
$\mathbf H_\alpha$ the resulting bundle. Local triviality at the gluing
points follows by extending $\alpha$ to a neighborhood of $A$, as
just justified. Define
\[
 d(\mathbf E,\mathbf F,\alpha)
 :=[\mathbf H_\alpha]-[r^*\mathbf F]\in K^0(D).
\]
This class lies in the kernel of $i_-^*$, because both bundles
restrict to $\mathbf F$ on $X_-$. The construction respects direct
sums and isomorphisms of triples. If $\alpha$ extends to an
isomorphism $\mathbf E\to\mathbf F$ on $X$, this isomorphism on $X_+$ and the
identity on $X_-$ give $\mathbf H_\alpha\cong r^*\mathbf F$;
thus an elementary triple has zero image.

A homotopy $(\alpha_t)_{t\in[0,1]}$ of isomorphisms over $A$ allows us to
glue the bundles over $D\times[0,1]$. The product is compact Hausdorff, and
extension from $A\times[0,1]$ is justified by the same argument.
Theorem~\ref{teo:homotopia-pullback-apendice} identifies the
restrictions of the glued bundle at the endpoints. Hence $d$ is unchanged
under the homotopy relation.

Let us also verify the composition relation. Let
$\alpha:\mathbf E|_A\to\mathbf F|_A$ and
$\beta:\mathbf F|_A\to\mathbf G|_A$. For
$0\leq t\leq\displaystyle\frac{\pi}{2}$ consider the isomorphisms
\[
 h_t(u,v):=
 \bigl(\cos t\,\alpha(u)+\sin t\,v,
       -\sin t\,\beta\alpha(u)+\cos t\,\beta(v)\bigr)
\]
from $(\mathbf E\oplus\mathbf F)|_A$ to
$(\mathbf F\oplus\mathbf G)|_A$. They are invertible because they are obtained
by composing $\alpha\oplus I$, the rotation of
$\mathbf F|_A\oplus\mathbf F|_A$ through angle $t$, and $I\oplus\beta$.
At $t=0$ we have $h_0=\alpha\oplus\beta$; at $t=\displaystyle\frac{\pi}{2}$,
composition with the global isomorphism
$(v,w)\mapsto(-w,v)$ from $\mathbf F\oplus\mathbf G$ to
$\mathbf G\oplus\mathbf F$ turns $h_{\displaystyle\frac{\pi}{2}}$ into
$(\beta\alpha)\oplus I_{\mathbf F}$. Gluing and homotopy give
\[
 \mathbf H_\alpha\oplus\mathbf H_\beta
 \cong\mathbf H_{\beta\alpha}\oplus r^*\mathbf F.
\]
Subtracting the classes of $r^*\mathbf F$ and $r^*\mathbf G$ yields
\[
 d(\mathbf E,\mathbf F,\alpha)+d(\mathbf F,\mathbf G,\beta)
 =d(\mathbf E,\mathbf G,\beta\alpha).
\]
Consequently, $d$ is well defined on the group of triples.

Let us prove its surjectivity onto the kernel of $i_-^*$. Let
$c\in K^0(D)$ with $i_-^*c=0$. By
Proposition~\ref{prop:igualdad-estable-k0}, we may write
$c=[\mathbf V]-[\underline{\mathbb C}^{\,N}_D]$. The equality
$i_-^*c=0$ allows us, after adding the same trivial bundle to both
sides, to choose a trivialization
\[
 \tau:\mathbf V|_{X_-}\xrightarrow{\;\cong\;}
 \underline{\mathbb C}^{\,N}_{X_-},
\]
where we have retained the letters $\mathbf V$ and $N$ for the stabilized
data. Take $\mathbf E:=\mathbf V|_{X_+}$,
$\mathbf F:=\underline{\mathbb C}^{\,N}_X$, and $\alpha:=\tau|_A$.
The bundle $\mathbf H_\alpha$ is isomorphic to $\mathbf V$: use
the identity on $X_+$ and $\tau^{-1}$ on $X_-$; these agree along the gluing locus.
Thus $d(\mathbf E,\mathbf F,\alpha)=c$.

Finally, suppose that $d(\mathbf E,\mathbf F,\alpha)=0$.
The same proposition provides an integer $N$ and an isomorphism
\[
 \mathbf H_\alpha\oplus\underline{\mathbb C}^{\,N}_D
 \cong r^*\mathbf F\oplus\underline{\mathbb C}^{\,N}_D.
\]
Its restrictions to $X_+$ and $X_-$ are isomorphisms
\[
 a:\mathbf E\oplus\underline{\mathbb C}^{\,N}_X
 \longrightarrow\mathbf F\oplus\underline{\mathbb C}^{\,N}_X,
 \qquad
 b:\mathbf F\oplus\underline{\mathbb C}^{\,N}_X
 \longrightarrow\mathbf F\oplus\underline{\mathbb C}^{\,N}_X
\]
satisfying $a|_A=b|_A\circ(\alpha\oplus I)$. Thus
$b^{-1}a$ extends $\alpha\oplus I$ to all of $X$. The original triple,
stabilized by an identity triple, is elementary and represents zero.
Every element of the group is represented by a triple, since addition is
direct sum and the inverse interchanges the bundles. This proves
injectivity of $d$.

The composition $(q^*)^{-1}d$ is the asserted isomorphism. The preceding verification
also shows that the relations can be realized by
stabilizations with elementary triples, isomorphisms, and homotopies: the
composition relation is precisely the rotation just
constructed. To compare two triples, apply the zero-kernel criterion
to the sum of the first with the inverse of the second.

All constructions are compatible with pullbacks. Indeed, a
map of pairs induces a map of their doubles commuting with the
retractions and quotients, and pullback commutes with gluing.
Uniqueness of the preimage under $q^*$ gives naturality. For a compact
Hausdorff space $\Lambda$, the same argument applies to
$(X\times\Lambda,A\times\Lambda)$ and its product with $[0,1]$; no
CW structure on $\Lambda$ is required.
\end{proof}

We turn to the noncompact case. Let $Y$ be a locally compact,
Hausdorff, noncompact space. Its one-point compactification is
\[
 Y^+:=Y\cup\{\infty\},
\]
where neighborhoods of $\infty$ are complements of compact
subsets of $Y$. If $Y$ is already compact, take the disjoint union with a
base point.

\begin{definition}[Compactly supported $K$--theory]
\label{def:k0-soporte-compacto}
\label{def:k-teoria-soporte-compacto}
Define
\[
 K_c^0(Y):=\widetilde K^0(Y^+).
\]
When $Y$ is the interior of a compact CW pair, and in particular when
$Y$ is the total space of a vector bundle over a closed manifold, a
class can equivalently be represented by a complex
$(\mathbf{E},\mathbf{F},a)$ in which $\mathbf{E}$ and $\mathbf{F}$ are bundles over $Y$ and
\[
 a\colon \mathbf{E}\longrightarrow \mathbf{F}
\]
is a bundle homomorphism invertible outside some compact set
$C\subseteq Y$. The
same relations as in Definition~\ref{def:triple-relativo-k} apply, after
enlarging the compact set when necessary.
\end{definition}

The word support refers to the region where the two-term complex
\[
 0\longrightarrow \mathbf{E}\xrightarrow{\ a\ }\mathbf{F}\longrightarrow0
\]
fails to be exact. No compactly supported section is chosen. If $Y$ is
the interior of a compact space $X$ with boundary $A=X\setminus Y$ and the
collapse map $X/A\to Y^+$ is a homeomorphism, then
\begin{equation}
 K_c^0(Y)\cong K^0(X,A).
 \label{eq:k-compacto-como-relativo}
\end{equation}

If $f\colon Y\to Z$ is proper, it extends continuously to
$f^+\colon Y^+\to Z^+$ by sending $\infty$ to $\infty$. Pullback defines
\[
 f^*\colon K_c^0(Z)\longrightarrow K_c^0(Y).
\]
Properness is necessary: without it, the preimage of the compact
support of a representative need not be compact.

\subsection{Suspension and the exact sequence}

For a pointed space $X$, its reduced suspension is
\[
 \Sigma X:=S^1\wedge X,
\]
where $\wedge$ denotes the smash product. We will need the odd group only
to write the exact sequence of a pair.

\begin{definition}[The group $K^1$]
\label{def:k1-suspension}
For a compact pointed space $X$, define
\[
 \widetilde K^1(X):=\widetilde K^0(\Sigma X).
\]
For an unpointed compact space, adjoin a disjoint base point and
set $K^1(X):=\widetilde K^1(X_+)$. For a locally compact space,
\[
 K_c^1(Y):=\widetilde K^1(Y^+).
\]
For a compact pair, set
\[
 K^1(X,A):=\widetilde K^1(X/A).
\]
\end{definition}

The definition by suspension agrees with the description that we
will use when applying periodicity. Indeed,
Proposition~\ref{prop:k1-producto-r-apendice} provides a natural
isomorphism
\begin{equation}
 K_c^1(Y)\cong K_c^0(\mathbb R\times Y).
 \label{eq:k1-como-producto-r}
\end{equation}
This is not a second definition: the pointed homeomorphism
$(\mathbb R\times Y)^+\cong S^1\wedge Y^+$ identifies the right-hand side
with the suspension appearing in Definition~\ref{def:k1-suspension}.

\begin{theorem}[Exactness for a pair]
\label{teo:sucesion-exacta-k-par}
Let $X$ be compact Hausdorff and let $A\subseteq X$ be closed. There exist connecting homomorphisms for
which the cyclic sequence
\[
 \begin{tikzcd}[column sep=large,row sep=large]
 K^0(X,A) \arrow[r] & K^0(X) \arrow[r] & K^0(A) \arrow[d,"\partial"] \\
 K^1(A) \arrow[u,"\partial"] & K^1(X) \arrow[l] & K^1(X,A) \arrow[l]
 \end{tikzcd}
\]
is exact. In particular,
\begin{equation}
 \operatorname{im}\bigl(K^0(X,A)\to K^0(X)\bigr)
 =\ker\bigl(K^0(X)\to K^0(A)\bigr).
 \label{eq:exactitud-k0-par-tramo}
\end{equation}
\end{theorem}

The proof belongs to algebraic topology: it identifies the connecting
maps with suspensions of quotient maps and uses Bott periodicity
to close the sequence. The statement for compact Hausdorff pairs with a closed subspace and
references for the proof are included in
Theorem~\ref{teo:exactitud-excision-k-apendice} of
Appendix~\ref{ap:topologia}. For the elementary construction of the segment
\eqref{eq:exactitud-k0-par-tramo}, see
\cite[Exercise~10.12(c), pp.~264--265]{BleeckerBoossIndex}; for the
compactly supported formulation, see
\cite[§10.4, pp.~266--269]{BleeckerBoossIndex}.

\section{External product}
\label{sec:producto-externo-k}

Let $X$ and $Y$ be compact, and let
$p_X\colon X\times Y\to X$ and $p_Y\colon X\times Y\to Y$ be the projections.
For bundles $\mathbf{E}\to X$ and $\mathbf{G}\to Y$ define
\[
 \mathbf{E}\boxtimes \mathbf{G}:=p_X^*\mathbf{E}\otimes p_Y^*\mathbf{G}.
\]
The fiber over $(x,y)$ is $\mathbf{E}_x\otimes \mathbf{G}_y$. Distributivity of tensor
product gives a bilinear map on the monoids of bundles and, by the
universal property of Grothendieck completion, an external product
\[
 \boxtimes\colon K^0(X)\otimes_{\mathbb Z}K^0(Y)
 \longrightarrow K^0(X\times Y).
\]
On representatives,
\begin{equation}
 \begin{aligned}
 \bigl([\mathbf{E}]-[\mathbf{F}]\bigr)\boxtimes\bigl([\mathbf{G}]-[\mathbf{H}]\bigr)
 ={}&[\mathbf{E}\boxtimes \mathbf{G}]+[\mathbf{F}\boxtimes \mathbf{H}]\\
 &-[\mathbf{E}\boxtimes \mathbf{H}]-[\mathbf{F}\boxtimes \mathbf{G}].
 \end{aligned}
 \label{eq:producto-externo-representantes}
\end{equation}

The formula proves directly that replacing either
representative by a stably equivalent one does not alter the result:
a summand added to both the positive and negative parts appears twice
with opposite signs. It also proves associativity. The interchange
isomorphism $\mathbf{E}_x\otimes \mathbf{G}_y\to \mathbf{G}_y\otimes \mathbf{E}_x$ gives commutativity after
interchanging the factors $X$ and $Y$.

\begin{proposition}[External product with compact support]
\label{prop:producto-externo-k}
Let $X$ and $Y$ be locally compact, Hausdorff, paracompact spaces.
There exists a natural bilinear
product
\[
 \boxtimes\colon K_c^0(X)\otimes_{\mathbb Z}K_c^0(Y)
 \longrightarrow K_c^0(X\times Y).
\]
If $a$ and $b$ admit representatives whose supports are contained in the
compact sets $C\subseteq X$ and $D\subseteq Y$, respectively, then
$a\boxtimes b$ admits a representative supported in
$C\times D$.
\end{proposition}

\begin{proof}
Existence, bilinearity, and naturality are the topological
properties of the compactly supported product stated in
Theorem~\ref{teo:exactitud-excision-k-apendice}. In the compactification
model, the reduced product on $X^+\times Y^+$ vanishes on
both axes and descends to the smash product
$X^+\wedge Y^+\cong(X\times Y)^+$.

We develop the support assertion in the complex model. Take two-term
representatives $a\colon\mathbf E\to\mathbf F$ and
$b\colon\mathbf G\to\mathbf H$, invertible outside $C$ and $D$, and fix
continuous Hermitian bundle metrics, which exist by paracompactness.
The even and odd parts of the graded product are
\[
 \mathbf V^+=(\mathbf E\boxtimes\mathbf G)
                    \oplus(\mathbf F\boxtimes\mathbf H),
 \qquad
 \mathbf V^-=(\mathbf F\boxtimes\mathbf G)
                    \oplus(\mathbf E\boxtimes\mathbf H).
\]
The product differential plus its adjoint, from the even part to the
odd part, is the homomorphism
\[
 D^+:=\begin{pmatrix}
       a\otimes I_{\mathbf G}&-I_{\mathbf F}\otimes b^*\\
       I_{\mathbf E}\otimes b&a^*\otimes I_{\mathbf H}
       \end{pmatrix}
       \colon\mathbf V^+\longrightarrow\mathbf V^-.
\]
Operators on distinct factors commute. Expanding the squares, the two
mixed terms cancel by the adjoint identities, and in each fiber we get
\[
 \begin{aligned}
 \|D^+(u,v)\|^2
 &=\|(a\otimes I_{\mathbf G})u\|^2
   +\|(I_{\mathbf E}\otimes b)u\|^2\\
 &\quad+\|(I_{\mathbf F}\otimes b^*)v\|^2
   +\|(a^*\otimes I_{\mathbf H})v\|^2.
 \end{aligned}
\]
If $x\notin C$, then $a_x$ and $a_x^*$ are invertible, and the first and
fourth norms show that $D^+_{(x,y)}$ is injective. If $y\notin D$, the
second and third norms give the same conclusion. In either case, the
even and odd parts have equal dimensions because at least one of the
pairs $\mathbf E_x,\mathbf F_x$ or $\mathbf G_y,\mathbf H_y$ has equal
ranks. Thus $D^+$ is an isomorphism outside $C\times D$. The product
representative $(\mathbf V^+,\mathbf V^-,D^+)$ is supported in that
compact set, as required.

\end{proof}

This construction agrees with that in
\cite[Theorem~10.16, pp.~267--268]{BleeckerBoossIndex}. In particular, the
Bott class on $\mathbb R^2$ can be multiplied externally by a
class on any locally compact space.

\section{The class of an elliptic symbol}
\label{sec:clase-simbolo-eliptico-k}

Let $(M,\mathbf{g})$ be a closed smooth Riemannian manifold. Denote by
$\pi\colon T^*M\to M$ the projection and by
\[
 B^*M:=\{(x,\xi)\in T^*M\mid \|\xi\|_{\mathbf{g}}\leq1\},
 \qquad
 S^*M:=\{(x,\xi)\in T^*M\mid \|\xi\|_{\mathbf{g}}=1\}
\]
the cotangent disk and sphere bundles. The radial map
\[
 \begin{aligned}
 \rho\colon T^*M&\longrightarrow\operatorname{int}(B^*M),\\
 (x,\xi)&\longmapsto
 \left(x,\frac{\xi}{\sqrt{1+\|\xi\|_{\mathbf{g}}^2}}\right)
 \end{aligned}
\]
is a diffeomorphism. Consequently,
\begin{equation}
 K_c^0(T^*M)\cong K^0(B^*M,S^*M).
 \label{eq:k-cotangente-relativo}
\end{equation}

Let $\mathbf{E},\mathbf{F}\to M$ be complex bundles and let
\[
 P\in\Psi^m_{\mathrm{cl}}(M;\mathbf{E},\mathbf{F})
\]
be a classical elliptic operator. By
Definition~\ref{def:simbolo-principal-pseudodiferencial-global},
\[
 \boldsymbol{\sigma}_m^\Psi(P)(x,t\xi)
 =t^m\boldsymbol{\sigma}_m^\Psi(P)(x,\xi),
 \qquad t>0,
\]
and ellipticity says that this bundle homomorphism is invertible when
$\xi\neq0$. In particular, its restriction to $S^*M$ is an isomorphism.

\begin{definition}[Symbol class]
\label{def:clase-simbolo-k}
Let $M$ be a closed smooth manifold, let
$\mathbf{E},\mathbf{F}\to M$ be finite-rank smooth complex vector
bundles, and let $P\in\Psi^m_{\mathrm{cl}}(M;\mathbf{E},\mathbf{F})$ be elliptic.
The \textbf{$K$--theory class of the symbol} of $P$ is
\begin{equation}
 [\boldsymbol{\sigma}(P)]
 :=
 \bigl[\pi^*\mathbf{E},\pi^*\mathbf{F},
 \boldsymbol{\sigma}_m^\Psi(P)|_{S^*M}\bigr]
 \in K^0(B^*M,S^*M)
 \cong K_c^0(T^*M).
 \label{eq:def-clase-simbolo-k}
\end{equation}
\end{definition}

The notation does not retain the order $m$: the class records only the isomorphism
on nonzero cotangent directions. For a differential operator,
the identity
\[
 \boldsymbol{\sigma}_m^\Psi(P)=i^m\boldsymbol{\sigma}_m(P)
\]
from \eqref{eq:conversion-simbolos-diferencial-pseudo} does not change the class. Indeed,
$i^m\in\mathbb C^\times$ can be joined to $1$ by a path in
$\mathbb C^\times$, and homotopy invariance of triples applies.

\begin{proposition}[Independence of choices]
\label{prop:independencia-clase-simbolo}
Let $M$ be a closed smooth manifold, let
$\mathbf{E},\mathbf{F}\to M$ be finite-rank smooth complex vector
bundles, and let $P\in\Psi^m_{\mathrm{cl}}(M;\mathbf{E},\mathbf{F})$ be a
classical elliptic operator.
The class \eqref{eq:def-clase-simbolo-k} is independent of the metric used to
define $B^*M$ and $S^*M$, after the canonical identification with
$K_c^0(T^*M)$. It is also unchanged when the symbol is replaced by another
elliptic symbol homotopic to it through elliptic symbols.
\end{proposition}

\begin{proof}
Let $\mathbf{g}_0$ and $\mathbf{g}_1$ be two metrics. For $\xi\neq0$,
define
\[
 q(x,\xi):=\frac{\|\xi\|_{\mathbf g_0}}{\|\xi\|_{\mathbf g_1}}>0,
 \qquad
 F(x,\xi):=(x,q(x,\xi)\xi).
\]
Extend $F$ by the identity on the zero section. Since $q$ is homogeneous
of degree zero and bounded above and below by positive constants, $F$
is a homeomorphism carrying
$(B^*_{\mathbf g_0}M,S^*_{\mathbf g_0}M)$ onto
$(B^*_{\mathbf g_1}M,S^*_{\mathbf g_1}M)$. On
$S^*_{\mathbf g_0}M$, homogeneity gives
\[
 \boldsymbol\sigma_m^\Psi(P)(F(x,\xi))
 =\|\xi\|_{\mathbf g_1}^{-m}
       \boldsymbol\sigma_m^\Psi(P)(x,\xi).
\]
The factor is positive and is joined to $1$ by a homotopy of boundary
isomorphisms:
\[
 \left[(1-t)+t\|\xi\|_{\mathbf g_1}^{-m}\right]
       \boldsymbol\sigma_m^\Psi(P)(x,\xi),\qquad0\leq t\leq1.
\]
This identifies the two triples after transporting the second by $F$.

We check that this transport induces the canonical identification with
$K_c^0(T^*M)$. If $\rho_i$ is the radial diffeomorphism used in
\eqref{eq:k-cotangente-relativo} for the metric $\mathbf g_i$, then
$F\rho_0=\rho_1F$, since
$\|q(x,\xi)\xi\|_{\mathbf g_1}=\|\xi\|_{\mathbf g_0}$.
Moreover, the maps
\[
 F_t(x,\xi):=(x,[(1-t)+tq(x,\xi)]\xi),\qquad0\leq t\leq1,
\]
extended over the zero section are homeomorphisms joining the identity
to $F$. Compactness of $M$ and uniform equivalence of the two norms
give constants $0<c\leq C$ such that, when
$(x,\xi_t)=F_t(x,\xi)$,
\[
 c\|\xi\|_{\mathbf g_0}
 \leq\|\xi_t\|_{\mathbf g_0}
 \leq C\|\xi\|_{\mathbf g_0}
 \qquad(x\in M,\ \xi\in T_x^*M,\ 0\leq t\leq1).
\]
By the lower bound, the preimage of a compact set under the homotopy is
contained in a fixed-radius disk bundle over $M$; it is closed and hence
compact. Thus the homotopy is proper and extends to the one-point
compactification while fixing the added point. Homotopy invariance
implies that $F$ induces the identity on $K_c^0(T^*M)$. The classes
obtained with the two metrics agree under the canonical identification.

If $a_t$ is a continuous family of invertible symbols over $S^*M$, the
triple over $(B^*M\times[0,1],S^*M\times[0,1])$ restricts at the endpoints to
the triples defined by $a_0$ and $a_1$. Homotopy invariance of
relative $K$--theory identifies their classes.
\end{proof}

Compact support can also be described using a symbol on
the entire cotangent bundle. The normalized
symbol
\[
 a(x,\xi)
 :=
 \|\xi\|_{\mathbf{g}}^{-m}\boldsymbol{\sigma}_m^\Psi(P)(x,\xi),
 \qquad \xi\neq0,
\]
is homogeneous of degree zero, but does not yet define a bundle homomorphism on the zero
section. Choose $\chi\in C^\infty([0,\infty))$ with $\chi=0$ for
$r\leq\varepsilon$ and $\chi=1$ for $r\geq2\varepsilon$, and set
\[
 \widetilde a(x,\xi)=
 \begin{cases}
 \chi(\|\xi\|_{\mathbf{g}})\|\xi\|_{\mathbf{g}}^{-m}\boldsymbol{\sigma}_m^\Psi(P)(x,\xi),
      &\xi\neq0,\\
 0,&\xi=0.
 \end{cases}
\]
Since $\widetilde a$ vanishes on a neighborhood of the zero section, it is smooth
on all of $T^*M$; outside the disk of radius $2\varepsilon$ it agrees with $a$ and
is invertible. The complex
$(\pi^*\mathbf{E},\pi^*\mathbf{F},\widetilde a)$ represents the symbol class. Two cutoffs or
two radii can be joined by a homotopy of cutoffs, retaining
invertibility outside a common disk.

\begin{theorem}[Realization by elliptic symbols]
\label{teo:realizacion-clases-simbolo}
Let $M$ be a closed smooth manifold. For every class
\[
 \kappa\in K_c^0(T^*M)
\]
there exist smooth complex bundles $\mathbf{E},\mathbf{F}\to M$ and an
elliptic operator
\[
 P\in\Psi^0_{\mathrm{cl}}(M;\mathbf{E},\mathbf{F})
\]
such that $[\boldsymbol{\sigma}(P)]=\kappa$.
\end{theorem}

\begin{proof}
Using \eqref{eq:k-cotangente-relativo}, represent $\kappa$ by a
triple $(\mathbf{V},\mathbf{W},\alpha)$ over $(B^*M,S^*M)$. The projection
$\pi\colon B^*M\to M$ and the zero section $\mathbf{s}\colon M\to B^*M$ are homotopy
inverses, since
\[
 H_t(x,\xi)=(x,t\xi),
 \qquad 0\leq t\leq1,
\]
contracts every disk to its center. Homotopy invariance of bundles
gives isomorphisms
\[
 \mathbf{V}\cong\pi^*\mathbf{E},
 \qquad
 \mathbf{W}\cong\pi^*\mathbf{F},
 \qquad
 \mathbf{E}:=\mathbf{s}^*\mathbf{V},
 \quad \mathbf{F}:=\mathbf{s}^*\mathbf{W}.
\]
The smoothing theorem for bundles over manifolds allows us to choose
smooth structures on $\mathbf{E}$ and $\mathbf{F}$ and to take the two preceding isomorphisms
smooth. This result is stated in Appendix~\ref{ap:topologia}.
Fix Hermitian metrics $\mathbf{h}_{\mathbf{E}}$ and $\mathbf{h}_{\mathbf{F}}$ and equip
$\operatorname{Hom}(\pi^*\mathbf{E},\pi^*\mathbf{F})$ with the induced operator norm. After
transporting $\alpha$, we obtain a continuous isomorphism
\[
 \alpha\colon\pi^*\mathbf{E}|_{S^*M}\longrightarrow\pi^*\mathbf{F}|_{S^*M}.
\]
If $S^*M=\varnothing$, then $M$ consists of finitely many points.
Take $P=0\colon\Gamma(\mathbf E)\to\Gamma(\mathbf F)$. This operator
has a smooth kernel and belongs to $\Psi^0_{\mathrm{cl}}$; the ellipticity
condition is vacuous, and its symbol triple is
$(\pi^*\mathbf E,\pi^*\mathbf F,\alpha)$, with the unique isomorphism
on the empty boundary. It therefore represents $\kappa$.

Suppose $S^*M\neq\varnothing$. Continuity of inversion and compactness
give
\[
 M_\alpha:=\max_{(x,\xi)\in S^*M}
       \|\alpha(x,\xi)^{-1}\|_{\operatorname{Hom}(\mathbf F_x,\mathbf E_x)}
       <\infty,
 \qquad\delta:=\frac{1}{2(1+M_\alpha)}>0.
\]
Choose a smooth approximation $\alpha_1$ such that
\[
 \sup_{(x,\xi)\in S^*M}
 \|\alpha_1(x,\xi)-\alpha(x,\xi)\|_{\operatorname{Hom}(\mathbf E_x,\mathbf F_x)}
 <\delta.
\]
For $\alpha_t:=\alpha+t(\alpha_1-\alpha)$, $0\leq t\leq1$, we have
\[
 \|\alpha(x,\xi)^{-1}
       (\alpha_t(x,\xi)-\alpha(x,\xi))\|_{\operatorname{End}(\mathbf E_x)}
 \leq M_\alpha\sup_{S^*M}\|\alpha_1-\alpha\|
 \leq\frac{M_\alpha}{2(1+M_\alpha)}<\frac{1}{2}.
\]
The Neumann series shows that $\alpha_t$ is invertible. This bound also
includes fibers of rank zero. Thus $\alpha$ and $\alpha_1$ define
homotopic triples; replace $\alpha$ by $\alpha_1$.

Extend $\alpha$ to $T^*M\setminus0_M$ by homogeneity of degree zero:
\[
 a_0(x,\xi)
 :=\alpha\left(x,\frac{\xi}{\|\xi\|_{\mathbf{g}}}\right).
\]
This section is smooth and invertible. Surjectivity in the exact
principal symbol sequence, Theorem~\ref{teo:sucesion-exacta-simbolo-principal-pseudo},
provides a classical operator
$P\in\Psi^0_{\mathrm{cl}}(M;\mathbf{E},\mathbf{F})$ with
$\boldsymbol{\sigma}_0^\Psi(P)=a_0$. Since $a_0$ is invertible in every nonzero direction,
$P$ is elliptic and its symbol triple is the triple representing
$\kappa$.
\end{proof}

The smoothing used in the proof is an approximation result for
continuous maps between bundles over a compact manifold. Its statement and
topological reference are included in Appendix~\ref{ap:topologia}. The
analytic step, quantization of $a_0$, was already proved in
Chapter~\ref{cap:operadores-pseudodiferenciales}.
The topological reduction allowing the two bundles of the triple to be
represented as pullbacks, and the subsequent realization, are also stated in
\cite[Exercise~11.4, pp.~279--280]{BleeckerBoossIndex}.

\section{Operations on symbol classes}
\label{sec:operaciones-clases-simbolos}

The rules of symbol calculus translate into exact identities in
$K_c^0(T^*M)$. We will prove them in the model by triples, keeping track of the
bundles and homomorphisms in each construction.

\begin{proposition}[Sum and stabilization]
\label{prop:suma-estabilizacion-simbolo-k}
Let $M$ be a closed smooth manifold and let
$\mathbf{E},\mathbf{F},\mathbf{E}',\mathbf{F}',\mathbf{G}\to M$ be finite-rank
smooth complex vector bundles.
Let
\[
 P\in\Psi^m_{\mathrm{cl}}(M;\mathbf{E},\mathbf{F}),
 \qquad
 Q\in\Psi^m_{\mathrm{cl}}(M;\mathbf{E}',\mathbf{F}')
\]
be elliptic operators of the same order. Then
\[
 [\boldsymbol{\sigma}(P\oplus Q)]=[\boldsymbol{\sigma}(P)]+[\boldsymbol{\sigma}(Q)].
\]
For an elliptic operator of arbitrary order and a bundle $\mathbf{G}\to M$,
stabilization by the identity triple satisfies
\[
 [\boldsymbol{\sigma}(P)]+[(\pi^*\mathbf{G},\pi^*\mathbf{G},I_{\pi^*\mathbf{G}})]=[\boldsymbol{\sigma}(P)].
\]
After reducing $P$ to order zero, this equality is realized by
direct sum with the operator $I_{\mathbf{G}}$.
\end{proposition}

\begin{proof}
The first equality is the definition of addition of triples:
\[
 (\pi^*\mathbf{E},\pi^*\mathbf{F},
   \boldsymbol{\sigma}_m^\Psi(P)|_{S^*M})
 \oplus
 (\pi^*\mathbf{E}',\pi^*\mathbf{F}',
   \boldsymbol{\sigma}_m^\Psi(Q)|_{S^*M})
 =
 (\pi^*(\mathbf{E}\oplus \mathbf{E}'),\pi^*(\mathbf{F}\oplus \mathbf{F}'),
 \boldsymbol{\sigma}_m^\Psi(P)|_{S^*M}
 \oplus\boldsymbol{\sigma}_m^\Psi(Q)|_{S^*M}).
\]
The triple $(\pi^*\mathbf{G},\pi^*\mathbf{G},I)$ extends as an isomorphism over all of
$B^*M$ and therefore represents zero. To justify the last assertion
at the operator level, fix a nonnegative covariant Laplacian $L_{\mathbf{E}}$.
Theorem~\ref{pseudo:teo-potencias-complejas-clasicas} shows that
\[
 \Lambda_{\mathbf{E}}^{-m}:=(I+L_{\mathbf{E}})^{-\frac{m}{2}}
 \in\Psi^{-m}_{\mathrm{cl}}(M;\mathbf{E},\mathbf{E})
\]
is elliptic and that its principal symbol on $S^*M$ is the identity. Thus
\[
 P^{(0)}:=P\Lambda_{\mathbf{E}}^{-m}
 \in\Psi^0_{\mathrm{cl}}(M;\mathbf{E},\mathbf{F})
\]
has the same symbol class as $P$, and
$P^{(0)}\oplus I_{\mathbf{G}}$ realizes the stated stabilization.
\end{proof}

\begin{proposition}[Composition]
\label{prop:composicion-simbolo-k}
Let $M$ be a closed smooth manifold and let
$\mathbf{E},\mathbf{F},\mathbf{G}\to M$ be finite-rank smooth complex vector
bundles.
Let
\[
 P\in\Psi^m_{\mathrm{cl}}(M;\mathbf{E},\mathbf{F}),
 \qquad
 Q\in\Psi^\ell_{\mathrm{cl}}(M;\mathbf{F},\mathbf{G})
\]
be elliptic operators. Then
\begin{equation}
 [\boldsymbol{\sigma}(QP)]=[\boldsymbol{\sigma}(P)]+[\boldsymbol{\sigma}(Q)].
 \label{eq:clase-simbolo-composicion}
\end{equation}
\end{proposition}

\begin{proof}
Theorem~\ref{teo:sucesion-exacta-simbolo-principal-pseudo} gives
\[
 \boldsymbol{\sigma}_{m+\ell}^\Psi(QP)
 =\boldsymbol{\sigma}_\ell^\Psi(Q)\circ\boldsymbol{\sigma}_m^\Psi(P).
\]
Applying the composition relation
\eqref{eq:relacion-composicion-triples} with
\[
 \alpha=\boldsymbol{\sigma}_m^\Psi(P)|_{S^*M},
 \qquad
 \beta=\boldsymbol{\sigma}_\ell^\Psi(Q)|_{S^*M},
\]
yields precisely \eqref{eq:clase-simbolo-composicion}.
\end{proof}

\begin{proposition}[Adjoint]
\label{prop:adjunto-simbolo-k}
Let $M$ be a closed smooth manifold and let
$\mathbf{E},\mathbf{F}\to M$ be finite-rank smooth complex vector
bundles; fix Hermitian metrics on them and a positive density on
$M$. If $P\in\Psi^m_{\mathrm{cl}}(M;\mathbf{E},\mathbf{F})$ is elliptic, then
\[
 [\boldsymbol{\sigma}(P^*)]=-[\boldsymbol{\sigma}(P)].
\]
\end{proposition}

\begin{proof}
Write
\[
 a:=\boldsymbol{\sigma}_m^\Psi(P)|_{S^*M}\colon\pi^*\mathbf{E}\longrightarrow\pi^*\mathbf{F}.
\]
By the symbolic adjoint rule,
$\boldsymbol{\sigma}_m^\Psi(P^*)|_{S^*M}=a^*$. The inverse of the class of the triple
$(\pi^*\mathbf{E},\pi^*\mathbf{F},a)$ is represented by
$(\pi^*\mathbf{F},\pi^*\mathbf{E},a^{-1})$. It suffices to compare $a^*$ with $a^{-1}$.

The endomorphism $a^*a$ of $\pi^*\mathbf{E}$ is positive definite. Matrix functional
calculus depends smoothly on the point of $S^*M$, so
\[
 h_t:=(a^*a)^{-t}a^*,
 \qquad 0\leq t\leq1,
\]
is a smooth homotopy through isomorphisms
$\pi^*\mathbf{F}\to\pi^*\mathbf{E}$. At the endpoints,
\[
 h_0=a^*,
 \qquad
 h_1=(a^*a)^{-1}a^*=a^{-1}.
\]
Homotopy invariance of triples completes the proof.
\end{proof}

\begin{proposition}[Homotopy and lower-order terms]
\label{prop:homotopia-orden-inferior-simbolo-k}
Let $M$ be a closed smooth manifold, let
$\mathbf{E},\mathbf{F}\to M$ be finite-rank smooth complex vector
bundles, and let $P_t\in\Psi^m_{\mathrm{cl}}(M;\mathbf{E},\mathbf{F})$, $0\leq t\leq1$, be a continuous
family whose principal symbols form a family of isomorphisms
over $S^*M$. Then
\[
 [\boldsymbol{\sigma}(P_0)]=[\boldsymbol{\sigma}(P_1)].
\]
In particular, if $R\in\Psi^{m-1}_{\mathrm{cl}}(M;\mathbf{E},\mathbf{F})$ and $P$ is
elliptic, then
\[
 [\boldsymbol{\sigma}(P+R)]=[\boldsymbol{\sigma}(P)].
\]
\end{proposition}

\begin{proof}
The first assertion is homotopy invariance of the symbol
triple. For the second,
\[
 \boldsymbol{\sigma}_m^\Psi(P+R)=\boldsymbol{\sigma}_m^\Psi(P),
\]
by exactness of the principal symbol sequence; in fact, the two triples are
identical.
\end{proof}

Combining realization with these identities explains why
$K_c^0(T^*M)$ is the natural domain of any additive invariant of
elliptic operators that is stable under homotopy.

\begin{lemma}[Relative smoothing of invertible homotopies]
\label{lem:suavizacion-relativa-homotopias-simbolo}
Let $M$ be a closed smooth manifold, let
$\mathbf{E},\mathbf{F}\longrightarrow M$ be finite-rank smooth bundles, and let
\[
a\colon S^*M\times[0,1]\longrightarrow
\operatorname{Iso}(\pi^*\mathbf{E},\pi^*\mathbf{F})
\]
be a continuous homotopy with smooth endpoints. There exists a smooth
homotopy $\widetilde a$ with the same endpoints, consisting entirely of
isomorphisms. Moreover, $a$ and $\widetilde a$ are joined, relative to
the endpoints, by a homotopy of isomorphisms.
\end{lemma}

Relative smoothing allows continuous homotopies to be used in proving
independence of the index from the representative, considered in
\cite[Exercise~12.13, pp.~291--292]{BleeckerBoossIndex}.

\begin{proof}
If $S^*M=\varnothing$, the conclusions are immediate. Otherwise,
take a smooth function $\sigma\colon[0,1]\to[0,1]$ equal to $0$ near
$0$ and to $1$ near $1$. Set
\[
 b(x,\xi,t):=a(x,\xi,\sigma(t)),\qquad
 H(r;x,\xi,t):=a(x,\xi,(1-r)t+r\sigma(t)).
\]
The homotopy $H$ consists of isomorphisms, fixes the endpoints $t=0,1$,
and joins $a$ to $b$. Since the endpoints of $a$ are smooth and $b$ is
constant in $t$ near them, $b$ is smooth on a neighborhood of the
endpoints.

Fix bundle metrics. Continuity of inversion and compactness give
\[
 M_0:=\max_{S^*M\times[0,1]}\|b^{-1}\|_{\operatorname{op}}<\infty,
 \qquad\mu:=\frac{1}{1+M_0}>0.
\]
The relative Whitney approximation theorem, applied to the section $b$
of the homomorphism bundle over $S^*M\times[0,1]$, produces a smooth
section $\widetilde a$, equal to $b$ near the endpoints, with
\[
 \|\widetilde a-b\|_{C^0(S^*M\times[0,1],
          \operatorname{Hom}(\pi^*\mathbf E,\pi^*\mathbf F))}
 <\frac{\mu}{2}.
\]
For every $0\leq\theta\leq1$,
\[
 (1-\theta)b+\theta\widetilde a
 =b[I+\theta b^{-1}(\widetilde a-b)],\qquad0\leq\theta\leq1,
\]
and
\[
 \|\theta b^{-1}(\widetilde a-b)\|_{C^0}
 \leq M_0\|\widetilde a-b\|_{C^0}
 \leq\frac{M_0}{2(1+M_0)}<\frac{1}{2}.
\]
The Neumann series shows that the entire segment consists of
isomorphisms. The segment fixes the endpoints. Concatenating it with
$H$ joins the original homotopy $a$ to $\widetilde a$ relative to the
endpoints and proves the assertions.
\end{proof}

\begin{theorem}[The analytic index on the symbol class]
\label{teo:indice-analitico-sobre-k}
Let $M$ be a closed smooth manifold. There exists a unique homomorphism
\begin{equation}
 \operatorname{ind}_{\mathrm a}\colon K_c^0(T^*M)\longrightarrow\mathbb Z
 \label{eq:indice-analitico-sobre-k}
\end{equation}
such that, for every classical elliptic pseudodifferential operator $P$,
\[
 \operatorname{ind}_{\mathrm a}([\boldsymbol{\sigma}(P)])=\operatorname{ind}P.
\]
\end{theorem}

\begin{proof}
By Theorem~\ref{teo:realizacion-clases-simbolo}, every class is the
symbol of some elliptic operator. For a smooth triple
$(\pi^*\mathbf{E},\pi^*\mathbf{F},a)$ over $(B^*M,S^*M)$ define
\[
 I(\mathbf{E},\mathbf{F},a):=\operatorname{ind}\operatorname{Op}(a),
\]
where $\operatorname{Op}(a)$ is any elliptic quantization of order
zero with principal symbol $a$. Two quantizations differ by an operator
of order $-1$, compact from $L^2(M,\mathbf E)$ to $L^2(M,\mathbf F)$ by
Rellich; Corollary~\ref{cor:indice-analitico-simbolo} shows that $I$ is
independent of the quantization. We must verify that this integer respects the
relations in the model by triples.

First consider two isomorphic triples, after any stabilizations
forming part of the isomorphism. Write them as
\[
(\pi^*\mathbf{E},\pi^*\mathbf{F},a),
\qquad
(\pi^*\mathbf{E}',\pi^*\mathbf{F}',a')
\]
over $(B^*M,S^*M)$. There exist isomorphisms
$U\colon\pi^*\mathbf{E}\to\pi^*\mathbf{E}'$ and
$V\colon\pi^*\mathbf{F}\to\pi^*\mathbf{F}'$ over $B^*M$ such that
$a'U=Va$ on $S^*M$. Although $U$ and $V$ may depend on the covariable,
radial contraction defines
\[
U_t(x,\xi)=U(x,t\xi),
\qquad
V_t(x,\xi)=V(x,t\xi),
\qquad 0\leq t\leq1.
\]
Therefore,
\[
a_t:=V_t^{-1}a'U_t
\]
is a homotopy of invertible symbols over $S^*M$ joining $a$ to
$V_0^{-1}a'U_0$. At the second endpoint, $U_0$ and $V_0$ are pullbacks of
continuous bundle isomorphisms over $M$.
Theorem~\ref{teo:suavizacion-haces-morfismos-apendice} allows us to approximate them
uniformly by smooth isomorphisms $\widehat U_0$ and $\widehat V_0$; an approximation smaller
than the minimum of their singular values preserves invertibility, and the
straight segments give homotopies of isomorphisms. Concatenate these segments
with $(a_t)$. Both endpoints of the resulting homotopy are now smooth,
so
Lemma~\ref{lem:suavizacion-relativa-homotopias-simbolo} allows us to replace it
by a smooth homotopy of symbols. A realization of the endpoint
$\widehat V_0^{-1}a'\widehat U_0$ is obtained from a realization of $a'$
by composition on the left and right with the bundle isomorphisms; these
compositions are invertible operators of order zero and do not change the
index. Thus an isomorphism of triples preserves the associated integer, even
when its initial isomorphisms depend on $\xi$.

A continuous homotopy of the boundary isomorphism can be smoothed relative
to its endpoints using
Lemma~\ref{lem:suavizacion-relativa-homotopias-simbolo} and quantized, by
surjectivity of the principal symbol, as a family of elliptic
operators: fix a single cover, a single partition of unity, and the
same cutoffs in the construction of
Theorem~\ref{teo:sucesion-exacta-simbolo-principal-pseudo}; the quantization
formula is linear in the symbol and preserves continuity in the
parameter.
Corollary~\ref{cor:indice-analitico-simbolo} gives constancy of its index.
After reducing all operators to order zero as in
Proposition~\ref{prop:suma-estabilizacion-simbolo-k}, direct sum of
triples corresponds to direct sum of operators.
Proposition~\ref{prop:calculo-indice-fredholm} proves additivity.

Suppose that the boundary isomorphism
$\alpha\colon\pi^*\mathbf{E}|_{S^*M}\to\pi^*\mathbf{F}|_{S^*M}$ extends to an
isomorphism $A\colon\pi^*\mathbf{E}\to\pi^*\mathbf{F}$ over $B^*M$. Radial contraction
shows that $A|_{S^*M}$ is homotopic, through isomorphisms, to the pullback
of $A|_{0_M}\colon \mathbf{E}\to \mathbf{F}$. This isomorphism may only be continuous.
Approximate it uniformly by a smooth isomorphism
$\widehat A_0\colon \mathbf{E}\to \mathbf{F}$; if the error is smaller than the minimum of the
singular values of $A|_{0_M}$, the segment between them remains
invertible. Concatenate this segment with radial contraction and apply
Lemma~\ref{lem:suavizacion-relativa-homotopias-simbolo}. The smooth endpoint
$\pi^*\widehat A_0$ is the constant symbol of multiplication by
$\widehat A_0$, whose index is zero. By homotopy invariance, any
realization of this triple has index zero. Thus the third relation in
Definition~\ref{def:triple-relativo-k} is respected.

Finally, the composition relation is respected because, for
elliptic $P\colon \mathbf{E}\to \mathbf{F}$ and $Q\colon \mathbf{F}\to \mathbf{G}$,
\[
 \operatorname{ind}(QP)
 =\operatorname{ind}Q+\operatorname{ind}P
\]
by Proposition~\ref{prop:calculo-indice-fredholm}, while their
symbols satisfy \eqref{eq:clase-simbolo-composicion}. Consequently,
the index descends to the quotient defined by triples and determines the
homomorphism \eqref{eq:indice-analitico-sobre-k}. Realization of all
classes also proves uniqueness.
\end{proof}

Bott periodicity and the Thom isomorphism will allow us to construct
the topological homomorphism on the same group.

\chapter{Bott periodicity, the Thom isomorphism, and the Euclidean index}
\label{cap:bott-thom-indice-euclidiano}

The symbol class constructed in
Chapter~\ref{cap:k-teoria-clase-simbolo} belongs to a compactly supported
group,
\[
 [\boldsymbol{\sigma}(P)]\in K_c^0(T^*M).
\]
Two operations are needed to turn it into an integer. Bott periodicity
identifies the $K$--theory of a space with that of its product with
$\mathbb R^2$, and the Thom isomorphism allows this identification to be made
fiber by fiber in a complex vector bundle. Both operations are topological;
the analytic work consists in proving that quantization of a class and the
Fredholm index behave in the same way.

We begin with the orientation conventions that determine the
signs. We then construct the Thom isomorphisms and the pushforward
maps associated with embeddings, calculate the indices of two model operators, and
prove the Euclidean formula. The final result will be
\[
 \operatorname{ind}_{\mathrm a}(P)
 =(-1)^n\alpha_n\bigl([\boldsymbol{\sigma}(P)]\bigr),
 \qquad P\in\operatorname{Ell}_c^0(\mathbb R^n),
\]
where $\alpha_n:K_c^0(\mathbb R^{2n})\longrightarrow\mathbb Z$ is the
iterated Bott isomorphism with the normalization specified
below. This is the form of
\cite[Theorem~11.5, pp.~281--282]{BleeckerBoossIndex} compatible with the
Fourier conventions of Chapter~\ref{cap:operadores-pseudodiferenciales}.

\begin{semblanzaHistorica}{Bott and Thom: transporting topological information}
The periodicity discovered by Raoul Bott revealed a deep repetition in the homotopy groups of the classical groups and gave $K$--theory an extraordinarily rigid structure. The Thom isomorphism allows this periodicity to be applied fiber by fiber in a vector bundle oriented in the appropriate sense. In the proof of the index theorem, these ideas push a localized class from the cotangent bundle of a manifold into Euclidean space, where its index can be calculated using model operators.
\end{semblanzaHistorica}

\section{The Bott class and orientation conventions}
\label{sec:clase-bott-convenciones}

Identify $\mathbb R^2$ with $\mathbb C$ by
\[
 (x,\xi)\longmapsto x+i\xi
\]
and orient $\mathbb R^2$ by the ordered basis
$(\partial_x,\partial_\xi)$. Its one-point compactification is
$(\mathbb R^2)^+\cong S^2$. Write
\[
 S^2=D_0\cup_{S^1}D_\infty
\]
as the union of two closed disks. If $m\in\mathbb Z$, denote by
$\mathbf{E}_m\longrightarrow S^2$ the line bundle obtained by identifying the
trivial fibers over the equator through
\[
 (z,v)_0\sim (z,z^m v)_\infty,
 \qquad z\in S^1,\quad v\in\mathbb C.
\]
The fiber of $\mathbf{E}_m$ over the point added at infinity is identified with
$\mathbb C$. Thus the difference of two rank-one bundles defines a
reduced class.

\begin{definition}[Bott class]
\label{def:clase-bott}
The \textbf{Bott class} used in this book is
\[
 b:=[\mathbf{E}_{-1}]-[\underline{\mathbb C}]
 \in \widetilde K^0(S^2)=K_c^0(\mathbb R^2).
\]
If $X$ is locally compact, write
\[
 b\boxtimes a\in K_c^0(\mathbb R^2\times X),
 \qquad a\in K_c^0(X),
\]
for the external product of
Proposition~\ref{prop:producto-externo-k}.
\end{definition}

The choice of $\mathbf{E}_{-1}$ rather than $\mathbf{E}_1$ is essential to the sign of the index.
Indeed, $\mathbf{E}_1\otimes \mathbf{E}_{-1}$ is trivial and the computation of the ring of the sphere
gives $b^2=0$; see
\cite[Theorem~2.2, pp.~41--50]{HatcherVBKT}. Since
$[\mathbf{E}_{-1}]=1+b$, we obtain
\[
[\mathbf{E}_1]=[\mathbf{E}_{-1}]^{-1}=(1+b)^{-1}=1-b.
\]
Thus the opposite class
$-b=[\mathbf{E}_1]-[\underline{\mathbb C}]$ is also a generator of
$K_c^0(\mathbb R^2)$, but interchanging these generators changes all the
odd-degree signs that will appear later.

It is useful to collect the conventions maintained throughout the argument:
\begin{equation}
 (x,\xi)\longmapsto x+i\xi,
 \qquad
 J(u,v)=(-v,u),
 \qquad
 b=[\mathbf{E}_{-1}]-[\underline{\mathbb C}],
 \qquad
 \lambda_{\mathbb C}=-b.
 \label{eq:resumen-convenciones-bott-thom}
\end{equation}
With this choice, the localized operator representing $b$ will have index
$-1$, whereas the normal oscillator representing
$\lambda_{\mathbb C}$ will have index $+1$. This difference is not a
subsequent sign correction; it is a consequence of the fact that the two operators
represent opposite generators.

The deep topological result that we will use is stated, with its
local compactness and naturality hypotheses, in
Theorem~\ref{teo:bott-periodicidad-apendice}. A complete proof can be
found in
\cite[Theorems~10.20 and~10.22, pp.~270--272]{BleeckerBoossIndex}.

\begin{theorem}[Complex Bott periodicity]
\label{teo:periodicidad-bott-k}
For every locally compact space $X$ there exists
a natural isomorphism
\[
 \alpha_X:K_c^0(\mathbb R^2\times X)
 \longrightarrow K_c^0(X)
\]
whose inverse is
\[
 \beta_X:K_c^0(X)\longrightarrow K_c^0(\mathbb R^2\times X),
 \qquad
 \beta_X(a)=b\boxtimes a.
\]
In particular,
\[
 \alpha_{\{*\}}(b)=1\in K^0(\{*\})=\mathbb Z.
\]
These isomorphisms are compatible with external products: if $Y$ is another
locally compact space, $u\in K_c^0(\mathbb R^2\times X)$ and
$v\in K_c^0(Y)$, then
\begin{equation}
 \alpha_{X\times Y}(u\boxtimes v)
 =\alpha_X(u)\boxtimes v.
 \label{eq:bott-compatibilidad-producto}
\end{equation}
\end{theorem}

We do not reproduce the topological proof of periodicity here. The consequences
we need follow immediately, and it is useful to record them.

\begin{corollary}[Iterated Bott isomorphism]
\label{cor:iteracion-bott}
For $n\in\mathbb N$ there exists an isomorphism
\[
 \alpha_{n,X}:K_c^0(\mathbb R^{2n}\times X)
 \longrightarrow K_c^0(X)
\]
obtained by applying $\alpha$ successively to the $n$ factors
$\mathbb R^2$. Its inverse is
\[
 a\longmapsto b^{\boxtimes n}\boxtimes a,
 \qquad
 b^{\boxtimes n}:=
 \underbrace{b\boxtimes\cdots\boxtimes b}_{n\text{ factors}}.
\]
In particular,
\[
 K_c^0(\mathbb R^{2n})=\mathbb Z\,b^{\boxtimes n},
 \qquad
 \alpha_n(b^{\boxtimes n})=1,
\]
where $\alpha_n:=\alpha_{n,\{*\}}$.
\end{corollary}

\begin{proof}
The assertion follows by induction. To pass from $n$ to $n+1$, use
\[
 \mathbb R^{2(n+1)}=\mathbb R^2\times\mathbb R^{2n}
\]
and apply first $\alpha_{\mathbb R^{2n}\times X}$ and then
$\alpha_{n,X}$. The formula for the inverse follows from
\eqref{eq:bott-compatibilidad-producto} and associativity of the external
product.
\end{proof}

Periodicity also determines the compactly supported $K$--theory groups
in odd dimensions. Using the natural isomorphism
\eqref{eq:k1-como-producto-r}, the definition by suspension in the preceding
chapter can be written as
\[
 K_c^1(X)\cong K_c^0(\mathbb R\times X).
\]
Periodicity then states that $K_c^{j+2}(X)\cong K_c^j(X)$. In particular,
\[
 K_c^0(\mathbb R^{2n+1})=0,
 \qquad
 K_c^1(\mathbb R^{2n+1})\cong\mathbb Z.
\]
Only even degree will be used explicitly, because $T^*\mathbb R^n$ has
real dimension $2n$.

\section{The Thom class}
\label{sec:clase-thom-k}

Let $\pi:\mathbf{V}\longrightarrow X$ be a complex vector bundle of rank
$r$ over a compact Hausdorff space $X$, equipped with a Hermitian
bundle metric $\mathbf{h}_{\mathbf{V}}$. The following construction is performed
fiber by fiber. If the base is merely locally compact, the same construction
defines a Thom orientation acting on compactly supported classes;
the orientation class itself may have noncompact support because
its support is the entire zero section.

For $v\in \mathbf{V}_x$, let
\[
 \varepsilon_v:\Lambda^k\mathbf{V}_x\longrightarrow\Lambda^{k+1}\mathbf{V}_x,
 \qquad
 \varepsilon_v(\omega)=v\wedge\omega,
\]
and let $\iota_v:=\varepsilon_v^*$ be its adjoint with respect to the Hermitian
inner product induced on the exterior algebra. The creation and
annihilation identities are
\begin{align}
 \varepsilon_v^2&=0,
 &\iota_v^2&=0,
 &\varepsilon_v\iota_v+\iota_v\varepsilon_v
 &=\|v\|^2I.
 \label{eq:creacion-aniquilacion-thom}
\end{align}
Define
\[
 c_{\mathbf{V}}(v):=\varepsilon_v-\iota_v:
 \Lambda^{\mathrm{ev}}\mathbf{V}_x\longrightarrow\Lambda^{\mathrm{odd}}\mathbf{V}_x.
\]
By \eqref{eq:creacion-aniquilacion-thom},
\begin{equation}
 c_{\mathbf{V}}(v)^*c_{\mathbf{V}}(v)=\|v\|^2I,
 \qquad
 c_{\mathbf{V}}(v)c_{\mathbf{V}}(v)^*=\|v\|^2I.
 \label{eq:clifford-thom-invertible}
\end{equation}
Consequently, $c_{\mathbf{V}}(v)$ is an isomorphism if $v\neq0$, with inverse
$\|v\|^{-2}c_{\mathbf{V}}(v)^*$.

\begin{definition}[Complex Thom class]
\label{def:clase-thom-compleja}
The \textbf{Thom class} of $\mathbf{V}$ is the class supported on the zero section
represented by
\[
 \lambda_{\mathbf{V}}
 :=[\pi^*\Lambda^{\mathrm{ev}}\mathbf{V},
     \pi^*\Lambda^{\mathrm{odd}}\mathbf{V},c_{\mathbf{V}}]
 \in K_c^0(\mathbf{V}).
\]
Equivalently, if $B(\mathbf{V})$ and $S(\mathbf{V})$ are the unit ball and sphere
bundles, then
\[
 \lambda_{\mathbf{V}}
 =\chi\bigl(
 \pi^*\Lambda^{\mathrm{ev}}\mathbf{V},
 \pi^*\Lambda^{\mathrm{odd}}\mathbf{V};
 c_{\mathbf{V}}|_{S(\mathbf{V})}
 \bigr)
 \in K^0(B(\mathbf{V}),S(\mathbf{V})).
\]
\end{definition}

Equality \eqref{eq:clifford-thom-invertible} proves all the analytic
and algebraic facts needed for the class to be defined. The
assertion that this class generates the $K$--theory of each fiber in a way
compatible with variation of the base is the topological content of the Thom
isomorphism. It is stated in
Theorem~\ref{teo:thom-k-apendice}; see also
\cite[Theorem~12.8, pp.~289--290]{BleeckerBoossIndex}.

\begin{theorem}[Thom isomorphism in $K$--theory]
\label{teo:isomorfismo-thom-k}
Let $\pi:\mathbf{V}\longrightarrow X$ be a finite-rank complex vector bundle over
a locally compact space. The map
\[
 \operatorname{Th}_{\mathbf{V}}:K_c^0(X)\longrightarrow K_c^0(\mathbf{V}),
 \qquad
 \operatorname{Th}_{\mathbf{V}}(a)=\pi^*a\cdot\lambda_{\mathbf{V}},
\]
is an isomorphism. If $f:Y\longrightarrow X$ is continuous and
$f^{-1}(\operatorname{supp}a)$ is compact, then
\begin{equation}
 (f_{\mathbf{V}})^*\operatorname{Th}_{\mathbf{V}}(a)
 =\operatorname{Th}_{f^*\mathbf{V}}(f^*a),
 \label{eq:naturalidad-thom}
\end{equation}
where $f_{\mathbf{V}}:f^*\mathbf{V}\longrightarrow \mathbf{V}$ is the induced map. If $\mathbf{V}$ and $\mathbf{W}$ are
complex bundles over $X$, then, under the natural identification of the total
bundle $\mathbf{V}\oplus \mathbf{W}$ with the bundle $\pi_{\mathbf{V}}^*\mathbf{W}$ over $\mathbf{V}$,
\begin{equation}
 \lambda_{\mathbf{V}\oplus \mathbf{W}}
 =\lambda_{\mathbf{V}}\boxtimes_X\lambda_{\mathbf{W}},
 \qquad
 \operatorname{Th}_{\mathbf{V}\oplus \mathbf{W}}
 =\operatorname{Th}_{\pi_{\mathbf{V}}^*\mathbf{W}}\circ\operatorname{Th}_{\mathbf{V}}.
 \label{eq:multiplicatividad-thom}
\end{equation}
\end{theorem}

For a noncompact base, the notation $\pi^*a\cdot\lambda_{\mathbf{V}}$ means that the
Thom complex is restricted to the preimage of a compact set containing
$\operatorname{supp}a$. Its support is contained in the zero section over
that compact set and therefore defines an element of $K_c^0(\mathbf{V})$.

Multiplicativity of the classes can be seen directly before using
the fact that the Thom maps are isomorphisms. Indeed,
\[
 \Lambda^*(\mathbf{V}\oplus \mathbf{W})
 \cong\Lambda^*\mathbf{V}\widehat\otimes\Lambda^*\mathbf{W}
\]
as $\mathbb Z_2$--graded bundles, and the corresponding odd operator is
\[
 c_{\mathbf{V}\oplus \mathbf{W}}(v,w)
 =c_{\mathbf{V}}(v)\widehat\otimes I
  +\Gamma_{\mathbf{V}}\widehat\otimes c_{\mathbf{W}}(w),
\]
where $\Gamma_{\mathbf{V}}$ equals $1$ in even degree and $-1$ in odd degree. The two
summands anticommute; hence the product of the adjoint with the operator is
\[
 \bigl(\|v\|^2+\|w\|^2\bigr)I.
\]
This is precisely the product complex representing the right-hand side
of \eqref{eq:multiplicatividad-thom}.

Let us check the relationship between the Bott and Thom conventions. On a trivial complex line we use the complex structure
\[
 J(u_1,u_2)=(-u_2,u_1).
\]
The map $c_{\mathbb C}(u_1+iu_2)$ is multiplication by
$u_1+iu_2$ from even to odd degree. With the gluing convention of
Definition~\ref{def:clase-bott}, this class is
\begin{equation}
 \lambda_{\mathbb C}=-b.
 \label{eq:thom-linea-menos-bott}
\end{equation}
Indeed, in trivializations over the interior and exterior disks, the transition
map of the preceding difference is $z$, whereas $b$ was defined
by $z^{-1}$. Consequently,
\begin{equation}
 \lambda_{\mathbb C^r}=(-1)^r b^{\boxtimes r}.
 \label{eq:thom-trivial-signo}
\end{equation}
This normalization will give the oscillator operator associated with the
Thom class index $1$ and the operator associated with $b$ index $-1$.

Not every real vector bundle has a canonical Thom class in
complex $K$--theory. A real bundle of even rank requires a
$\operatorname{Spin}^c$ orientation; once it has been chosen, the spinor bundles
$\mathbf{S}^+$ and $\mathbf{S}^-$ and Clifford multiplication produce a class
\[
 [\pi^*\mathbf{S}^+,\pi^*\mathbf{S}^-,c(v)]\in K_c^0(\mathbf{V}).
\]
A complex bundle has a canonical $\operatorname{Spin}^c$ orientation, and the
resulting class agrees with the exterior-algebra class above. We will not apply the
Thom isomorphism directly to an arbitrary real bundle. The embedding
used for the index theorem will instead involve the sum of two
copies of a real bundle; this sum has a canonical complex structure and does not
require an additional choice of $\operatorname{Spin}^c$ structure.

\section{Embeddings and extension by zero}
\label{sec:pushforward-k-encajes}

Let $U\subseteq Z$ be an open subset of a locally compact space. We will use the
description by supports
\begin{equation}
 K_c^0(U)
 =\varinjlim_{K\Subset U}\mathcal K_K^0(U),
 \label{eq:k-compacto-colimite-soportes}
\end{equation}
with the notation by complexes invertible outside $K$ introduced in
Appendix~\ref{ap:topologia}. For every compact set $K\Subset U$, excision
states that restriction to $U$ and extension to $Z$ are inverse
isomorphisms. We will use the second map, written
\begin{equation}
 \mathcal K_K^0(U)
 \xrightarrow{\ \cong\ }
 \mathcal K_K^0(Z).
 \label{eq:excision-extension-cero-k}
\end{equation}
These identifications commute with the morphisms induced by inclusions
of compact sets. Passing to the direct limit defines extension by zero
\begin{equation}
 e_U:K_c^0(U)\longrightarrow K_c^0(Z).
 \label{eq:extension-cero-k}
\end{equation}
If a class is represented by a complex elementary outside a
compact set $K\Subset U$, its image is the relative class determined by the
same support $K$ through \eqref{eq:excision-extension-cero-k}. This
description avoids requiring the individual bundles of the representative to extend
to all of $Z$. Compatibility of excision with inclusions
gives
\[
 e_V\circ e_U=e_U
\]
for open subsets $U\subseteq V\subseteq Z$, where on the left-hand side
the first extension is from $U$ to $V$ and the second from $V$ to $Z$.

Now consider a proper smooth embedding
\[
 j:X\hookrightarrow Y
\]
between manifolds without boundary. Choose metrics and use their musical
isomorphisms to identify cotangent bundles with tangent bundles. The
differential induces an embedding
\[
 Tj:TX\hookrightarrow TY.
\]
Let $\nu\longrightarrow X$ be the normal bundle of $j$, obtained from the
chosen metric on $Y$. Identify $X$ with its image and write
$T_xY=T_xX\oplus\nu_x$. Fix a connection $\nabla^Y$ on $TY$.
The tangential projection of $\nabla^Y$ induces a connection $\nabla^X$
on $TX$, and its normal component defines the tensor
\[
 B_x(w,v):=\bigl(\nabla^Y_w V\bigr)^\perp,
 \qquad w,v\in T_xX,
\]
where $V$ is any vector field tangent to $X$ with $V(x)=v$. The expression
is independent of $V$: writing a field vanishing at $x$ in a tangent
frame, the terms containing derivatives of its coefficients
are tangent and disappear upon projection onto $\nu_x$. If $\nabla^Y$
is the Levi--Civita connection, $B$ is the second fundamental form.

Proposition~\ref{prop:descomposicion-horizontal-vertical-conexion},
applied to these two connections, gives decompositions
\[
 T_v(TY)\cong T_xY\oplus T_xY,
 \qquad
 T_v(TX)\cong T_xX\oplus T_xX,
 \qquad (x,v)\in TX.
\]
The tangent inclusion is generally not the direct inclusion of the two
summands. To calculate it, represent a vector in $T_v(TX)$ by a
curve $v(t)\in T_{x(t)}X$ and set
$w=x'(0)$, $z=\nabla^X_{\partial t}v(0)$. Decomposition of the
connection along the curve gives
$\nabla^Y_{\partial t}v(0)=z+B_x(w,v)$. Thus, in horizontal--vertical
coordinates, the differential of $Tj$ is
\[
 T_v(Tj)(w,z)=\bigl(w,z+B_x(w,v)\bigr).
\]
Now consider the linear map
\[
 F_v:T_xY\oplus T_xY\longrightarrow\nu_x\oplus\nu_x,
 \qquad
 F_v(a,b):=\bigl(a^\perp,b^\perp-B_x(a^\top,v)\bigr).
\]
It is surjective, since $F_v(u_1,u_2)=(u_1,u_2)$ for
$u_1,u_2\in\nu_x$. If $F_v(a,b)=0$, then $a\in T_xX$ and
$b=b^\top+B_x(a,v)$, so
$(a,b)=T_v(Tj)(a,b^\top)$. Conversely, the formula for $T_v(Tj)$
shows that its entire image lies in the kernel of $F_v$. Thus $F_v$
induces an isomorphism from the quotient by the image of $T_v(Tj)$ to
$\nu_x\oplus\nu_x$.

These maps depend smoothly on $(x,v)$, and their inverse in the
quotient is the class of $(u_1,u_2)$. Identifying the normal bundle with this
quotient by means of an orthogonal complement yields
\begin{equation}
 N(Tj)\cong\pi_{TX}^*(\nu\oplus\nu).
 \label{eq:normal-tx-en-ty}
\end{equation}
The identification depends on the connection. If two connections are joined by
their affine segment, the induced connections, tensors $B$, and
maps $F_v$ depend smoothly on the parameter, and each $F_v$ induces the
preceding isomorphism. This gives a homotopy of identifications.
On the right-hand side of \eqref{eq:normal-tx-en-ty}, fix the complex
structure
\begin{equation}
 J(u_1,u_2)=(-u_2,u_1).
 \label{eq:estructura-compleja-normal-tx}
\end{equation}
The first summand corresponds to the horizontal direction and the second to the
vertical direction. Fix this structure $J$ and a tubular neighborhood
$\Phi_J$ of $Tj$; the complex orientation, Thom isomorphism, and
extension by zero are always taken with respect to the same pair $(J,\Phi_J)$.
Thus, even if $\nu$ is a real bundle without a
$\operatorname{Spin}^c$ orientation, the normal bundle of $TX$ in $TY$ is a complex bundle.

A tubular neighborhood provides a diffeomorphism
\[
 \Phi_J:\mathcal U\subseteq N(Tj)\longrightarrow\mathcal V\subseteq TY
\]
between neighborhoods of the zero sections. By excision, every compactly
supported class in $N(Tj)$ can be represented inside $\mathcal U$ after
a radial dilation. We then compose the Thom isomorphism with
$\Phi_J$ and extension by zero from $\mathcal V$.

\begin{definition}[Pushforward of an embedding]
\label{def:pushforward-k-encaje}
The \textbf{pushforward} associated with $j$ is
\begin{equation}
 j_!:K_c^0(T^*X)\longrightarrow K_c^0(T^*Y),
 \qquad
 j_!:=e_{\mathcal V}\circ(\Phi_J^{-1})^*\circ
 \operatorname{Th}_{N(Tj)},
 \label{eq:def-pushforward-k-encaje}
\end{equation}
where the metrics have been used to identify tangent and cotangent bundles and
$N(Tj)$ carries the complex structure
\eqref{eq:estructura-compleja-normal-tx}.
\end{definition}

Independence of the metrics, connection, and tubular neighborhood follows
from homotopy invariance, naturality of Thom, and isotopy of tubular
neighborhoods. Full functoriality is the topological result
recorded in Theorem~\ref{teo:mapa-directo-k-apendice}; in the preceding notation
it states
\begin{equation}
 (k\circ j)_!=k_!\circ j_!
 \label{eq:functorialidad-pushforward-k}
\end{equation}
for composable proper embeddings $X\xrightarrow{j}Y\xrightarrow{k}Z$.
The construction in the form used here is described in
\cite[pp.~289--294]{BleeckerBoossIndex}; independence of choices and
functoriality are proved in
\cite[Appendix~C, Theorem~C.8; Chapter~III,
\S\S12--13]{LawsonMichelsohn1989}.

For the inclusion $i:\{0\}\hookrightarrow\mathbb R^N$, the normal bundle of
$T\{0\}$ in $T\mathbb R^N$ is
\[
 \mathbb R^N\oplus\mathbb R^N\cong\mathbb C^N
\]
with the complex structure of \eqref{eq:estructura-compleja-normal-tx}.
By \eqref{eq:thom-trivial-signo},
\begin{equation}
 i_!(1)=(-1)^N b^{\boxtimes N}.
 \label{eq:pushforward-punto-bott}
\end{equation}
Since $\alpha_N(b^{\boxtimes N})=1$, we obtain
\begin{equation}
 i_!^{-1}=(-1)^N\alpha_N:
 K_c^0(T\mathbb R^N)=K_c^0(\mathbb R^{2N})
 \longrightarrow\mathbb Z.
 \label{eq:inversa-pushforward-punto}
\end{equation}
This formula explains in advance the sign that will appear in the
Euclidean index.

\section{Euclidean operators equal to the identity at infinity}
\label{sec:ellc-euclidiano}

The phrase ``at infinity'' will refer to a strong support condition.
After stabilizing the bundles
by trivial summands, all operators will act on functions with
values in a single space $\mathbb C^q$.

\begin{definition}[The class $\operatorname{Ell}_c^0(\mathbb R^n)$]
\label{def:ellc-euclidiano}
An operator $P$ belongs to
$\operatorname{Ell}_c^0(\mathbb R^n;\mathbb C^q)$ if it is a
classical, properly supported, elliptic pseudodifferential operator of order zero,
and there exists a compact set $K\subset\mathbb R^n$ such that
\[
 P=I+Q,
 \qquad
 \operatorname{supp}(\mathbf{K}_Q)\subseteq K\times K,
\]
where $\mathbf{K}_Q$ is the Schwartz kernel of $Q$. Direct sum with
identity operators is allowed, and we write simply
$\operatorname{Ell}_c^0(\mathbb R^n)$ for the stable union of these
classes.
\end{definition}

The condition on the Schwartz kernel is equivalent to the two statements
\[
 (P-I)u=0,
 \qquad
 (P^*-I)u=0
\]
for every $u\in C_c^\infty(\mathbb R^n;\mathbb C^q)$ whose support does not meet
$K$. One implication is immediate. For the converse, the first equality
shows, in the distributional sense, that the kernel of $P-I$ vanishes when
the second variable lies outside $K$; the equality for the adjoint shows
that it vanishes when the first variable lies outside $K$. It also follows
that the principal symbol is the identity if $x\notin K$.

Equality of the principal symbol with $I$ outside a compact set alone is a
weaker condition: it does not assert that the operator is exactly the identity
there. To compare realizations of this symbol, first carry out the
quantization on $S^n=\mathbb R^n\cup\{\infty\}$. A matrix-valued
principal symbol $a_0$ equal to $I$ outside a compact set
$K\subset\mathbb R^n$ extends by $I$ to the cotangent bundle of $S^n$.
Let $I+A\in\Psi^0_{\mathrm{cl}}(S^n;\mathbb C^q)$ be a quantization
of this extension. Choose $\chi\in C_c^\infty(\mathbb R^n)$ equal to
$1$ on a neighborhood of $K$, extend it by zero to $S^n$, and
consider
\[
 I+\chi A\chi.
\]
Its principal symbol is $I+\chi^2(a_0-I)=a_0$. Consequently,
$A-\chi A\chi\in\Psi^{-1}_{\mathrm{cl}}(S^n)$ and it is compact on
$L^2(S^n)$: it maps $L^2(S^n)$ continuously into $H^1(S^n)$, whose
inclusion in $L^2(S^n)$ is compact. By
Theorem~\ref{pseudo:teo-estabilidad-indice}, the two operators on
$S^n$ have the same index. Moreover, the kernel of $\chi A\chi$ is
supported in $\operatorname{supp}\chi\times\operatorname{supp}\chi$,
so its Euclidean restriction satisfies the strong condition of
Definition~\ref{def:ellc-euclidiano}. This construction allows
realizations with strong support to be used in the stable theory described in
\cite[p.~275]{BleeckerBoossIndex}. The compactness comparison has been
made on $S^n$; equality of principal symbols alone does not
allow it to be applied to an arbitrary quantization on $\mathbb R^n$.

In this compactification, if
$P\in\operatorname{Ell}_c^0(\mathbb R^n)$, equality with the identity near infinity
allows it to be extended by $I$ to an elliptic pseudodifferential operator
\[
 \overline P\in\Psi^0_{\mathrm{cl}}(S^n).
\]
Choose $K$ larger if necessary so that it is a compact neighborhood of
the kernel support. The local order-zero mapping theorem shows that $Q$
is bounded on $L^2(\mathbb R^n;\mathbb C^q)$: insert compact cutoffs equal
to one on $K$ and use the Euclidean Sobolev mapping theorem. The support
condition implies $Q=1_KQ1_K$, first distributionally on test functions
and then on $L^2(\mathbb R^n;\mathbb C^q)$ by boundedness. Hence, relative to
\[
 L^2(\mathbb R^n;\mathbb C^q)
 =L^2(K;\mathbb C^q)\oplus L^2(K^c;\mathbb C^q),
 \qquad P=(I+Q)_K\oplus I.
\]
The extension $\overline P$ has the same block decomposition on $S^n$.
In stereographic coordinates its volume measure is $w(x)\,dx$ with $w$
smooth and positive. On the compact set $K$, $w$ is bounded above and
below by positive constants. Thus the two norms on the active block
$L^2(K;\mathbb C^q)$ are equivalent and the operator on this block is
literally the same. Fredholmness, including closedness of the range,
and the index are invariant under this equivalence of norms. Since
$\overline P$ is Fredholm by Theorem~\ref{pseudo:teo-fredholm-eliptico},
$P$ is Fredholm on $L^2(\mathbb R^n;\mathbb C^q)$ and
\begin{equation}
 \operatorname{ind}P=\operatorname{ind}\overline P.
 \label{eq:indice-compactificacion-euclidiana}
\end{equation}
This argument does not identify adjoints computed with different volume
measures.

The principal symbol of $P$ gives, by
Definition~\ref{def:clase-simbolo-k}, a class
\[
 [\boldsymbol{\sigma}(P)]\in K_c^0(T^*\mathbb R^n)
 =K_c^0(\mathbb R^{2n}).
\]
Its support is compact because the symbol is the identity outside $K$ in the
base variable and is invertible in every nonzero cotangent direction. In the
relative model it is represented by
\[
 [\underline{\mathbb C}^q,
   \underline{\mathbb C}^q,
   \boldsymbol{\sigma}_0^\Psi(P)|_{S^*\mathbb R^n}].
\]

\begin{proposition}[Euclidean realization and descent of the index]
\label{prop:indice-analitico-euclidiano-k}
Every class $a\in K_c^0(T^*\mathbb R^n)$ is the symbol class of some
$P\in\operatorname{Ell}_c^0(\mathbb R^n)$ after stabilization.
Moreover, there exists a unique homomorphism
\begin{equation}
 I_n:K_c^0(T^*\mathbb R^n)\longrightarrow\mathbb Z
 \label{eq:hom-indice-euclidiano}
\end{equation}
such that
\[
 I_n([\boldsymbol{\sigma}(P)])=\operatorname{ind}P.
\]
\end{proposition}

\begin{proof}
Represent $a$ by a compactly supported triple
\[
 (\mathbf{E},\mathbf{F},\theta)
\]
over $T^*\mathbb R^n$. Since
$T^*\mathbb R^n\cong\mathbb R^{2n}$ is contractible, after adding an
elementary triple we may assume that $\mathbf{E}$ and $\mathbf{F}$ are trivial bundles of
the same rank.

Extend $a$ by zero to the cotangent bundle of
$S^n=\mathbb R^n\cup\{\infty\}$. In the relative model, choose a
representative over $(B^*S^n,S^*S^n)$ elementary over $T^*U$ for
a neighborhood $U$ of $\infty$. The relative proof of
Theorem~\ref{teo:realizacion-clases-simbolo} allows us to smooth it without
changing that neighborhood. After stabilizing and fixing the
elementary trivialization there, we obtain an elliptic symbol of order zero
that is exactly the identity over $T^*U$.

Let us construct the operator representing this class. Let $a_0(x,\xi)$ be the resulting homogeneous
symbol, equal to $I$ when $x$ is near $\infty$.
Theorem~\ref{teo:sucesion-exacta-simbolo-principal-pseudo} gives a
quantization $A$ with principal symbol $a_0$. Choose
$\chi\in C_c^\infty(\mathbb R^n)$ equal to $1$ on the projection to the base
of $\operatorname{supp}(a_0-I)$ and define
\[
 P:=I+\chi(A-I)\chi.
\]
Its principal symbol is
\[
 I+\chi(x)^2(a_0(x,\xi)-I)=a_0(x,\xi),
\]
because $a_0-I$ vanishes where $\chi\neq1$. Without imposing an
additional condition on quantization, the two multiplications by $\chi$ show
directly that the Schwartz kernel of $P-I$ is contained in
$K\times K$, with $K=\operatorname{supp}\chi$. Therefore,
$P\in\operatorname{Ell}_c^0(\mathbb R^n)$ and
$[\boldsymbol{\sigma}(P)]=a$.

To prove that the index descends to $K$--theory, all the relative
information of the triple must be retained; an isomorphism of triples may depend
on $\xi$ and does not automatically become a conjugation of operators. Let
$a_0$ and $a_1$ be two representatives of the same class. After adding elementary
triples, the relations of the relative model provide a parametrized
triple over
\[
 B^*S^n\times[0,1]
\]
whose endpoint restrictions are the given representatives and which is
elementary over $B^*U\times[0,1]$ for a common neighborhood $U$ of $\infty$.
Relative smoothing keeps this neighborhood fixed. Radial contraction in
the fibers of $B^*S^n$ turns the stabilized isomorphisms into
bundle isomorphisms over the base, exactly as in the proof of
Theorem~\ref{teo:indice-analitico-sobre-k}. Homogeneous extension in the
covariable produces a family of elliptic symbols of order zero equal to
the identity over $T^*U$.

Using a single cover, a single partition of unity, and a single cutoff
$\chi$, quantization is continuous in the topology of symbols of order
zero and, by Calderón--Vaillancourt, continuous in the norm of
\[
 \mathcal L\bigl(L^2(\mathbb R^n;\mathbb C^q)\bigr).
\]
Each operator extends to an elliptic operator on $S^n$. Stability
of the Fredholm index gives
\[
 \operatorname{ind}P_0=\operatorname{ind}P_1.
\]
If the triple is elementary, radial contraction takes it to a bundle
isomorphism over the base; its quantization is an isomorphism of order zero plus
an operator of order $-1$ and has index zero. Finally, direct sum
of triples is quantized by direct sum of operators, and
\[
 \operatorname{ind}(P\oplus Q)
 =\operatorname{ind}P+\operatorname{ind}Q.
\]
Thus the index respects the stabilization, homotopy,
isomorphism, and addition relations defining $K_c^0(T^*\mathbb R^n)$. This determines
$I_n$, and realization of all classes proves its uniqueness.
\end{proof}

The same construction works for families. If $\Lambda$ is compact
Hausdorff and
$P_\lambda$ is a family continuous in the topology of classical symbols of
order zero, with a common compact support, then the family is
norm-continuous and its principal symbols depend jointly on
$(x,\xi,\lambda)$. Under these hypotheses, the symbols define
\begin{equation}
 [\boldsymbol{\sigma}(P_\Lambda)]
 \in K_c^0(T^*\mathbb R^n\times\Lambda).
 \label{eq:clase-familia-simbolos-euclidiana}
\end{equation}
Realization of every class and descent of the index bundle will be proved in
Proposition~\ref{prop:realizacion-descenso-familias-euclidianas}. The
variable $\lambda$ is not incorporated into the cotangent variables: it continuously
parametrizes the bundles and the symbol isomorphism.

\section{Two model operators}
\label{sec:operadores-modelo-bott-thom}

Periodicity determines the group
$K_c^0(T^*\mathbb R^n)$, but does not determine the analytic homomorphism
$I_n$ until its value on a generator is fixed. We first calculate that value
for the class $b$ using a localized Toeplitz operator. We then
calculate the oscillator operator representing the opposite class $-b$, which
will be the analytic model for the Thom isomorphism.

\subsection{The localized operator representing \texorpdfstring{$b$}{b}}

Let $S^1\subset\mathbb C$ and let
\[
 e_k(z)=z^k,
 \qquad k\in\mathbb Z,
\]
be the orthonormal Fourier basis of $L^2(S^1)$. Define
\[
 H_+:=\overline{\operatorname{span}}\{e_k\mid k\geq0\},
 \qquad
 H_-:=\overline{\operatorname{span}}\{e_k\mid k<0\},
\]
and denote by $\Pi_+$ and $\Pi_-$ the corresponding orthogonal
projections. The projection $\Pi_+$ is pseudodifferential of order zero:
modulo the finite-dimensional projection onto the constants,
\[
 \Pi_+=\frac{1}{2}\bigl(I+\operatorname{sgn}D\bigr),
 \qquad
 D=-i\frac{d}{d\theta},
\]
and its principal symbol equals $1$ on $\xi>0$ and $0$ on $\xi<0$.

Let $g:S^1\longrightarrow S^1$ be a smooth function of degree one equal to
$1$ on a neighborhood of a point $p\in S^1$. It can be obtained by deforming
$z\mapsto z$ through functions of degree one and concentrating all
angular variation in an arc disjoint from $p$. Let $M_g$ be multiplication
by $g$ and consider
\begin{equation}
 A_g:=\Pi_+M_g\Pi_++\Pi_-:
 L^2(S^1)=H_+\oplus H_-\longrightarrow L^2(S^1).
 \label{eq:toeplitz-bott}
\end{equation}
Its principal symbol is
\begin{equation}
 \boldsymbol{\sigma}_0^\Psi(A_g)(z,\xi)=
 \begin{cases}
 g(z),&\xi>0,\\
 1,&\xi<0.
 \end{cases}
 \label{eq:simbolo-toeplitz-bott}
\end{equation}
It is invertible and equals $1$ over $T^*U\setminus0$ for some neighborhood $U$ of
$p$. Thus $A_g$ is elliptic.

\begin{lemma}[Index and class of the Bott model]
\label{lem:modelo-bott-indice}
There exists $T_b\in\operatorname{Ell}_c^0(\mathbb R)$ such that
\[
 [\boldsymbol{\sigma}(T_b)]=b,
 \qquad
 \operatorname{ind}T_b=-1.
\]
\end{lemma}

\begin{proof}
We begin with the index of \eqref{eq:toeplitz-bott}. A smooth homotopy
$g_t$ between $g_0(z)=z$ and $g_1=g$ gives, through
\eqref{eq:toeplitz-bott}, a norm-continuous family of elliptic
operators. Its index is constant. For $g_0(z)=z$ we have
\[
 A_{g_0}e_k=
 \begin{cases}
 e_{k+1},&k\geq0,\\
 e_k,&k<0.
 \end{cases}
\]
Thus
\[
 \ker A_{g_0}=\{0\},
 \qquad
 \operatorname{Ran}A_{g_0}=\{e_0\}^{\perp},
\]
and consequently
\begin{equation}
 \operatorname{ind}A_g=\operatorname{ind}A_{g_0}=-1.
 \label{eq:indice-toeplitz-generador}
\end{equation}

Choose $\chi\in C^\infty(S^1)$ with support disjoint from $p$, equal to $1$
on the support of $g-1$, and define
\[
 \widetilde T_b:=I+\chi(A_g-I)\chi.
\]
By \eqref{eq:simbolo-toeplitz-bott}, $\widetilde T_b$ and $A_g$ have the
same principal symbol. Their difference has order $-1$ on the closed
manifold $S^1$ and is therefore compact on $L^2(S^1)$. By
Theorem~\ref{pseudo:teo-estabilidad-indice},
\[
 \operatorname{ind}\widetilde T_b=-1.
\]
Moreover, the kernel of $\widetilde T_b-I$ is contained in
$\operatorname{supp}\chi\times\operatorname{supp}\chi$. Removing the point
$p$ and identifying $S^1\setminus\{p\}$ with $\mathbb R$ gives an
operator $T_b\in\operatorname{Ell}_c^0(\mathbb R)$ with the same index.

It remains to identify its class. Compactify the phase plane
$T^*\mathbb R\cong\mathbb R_x\times\mathbb R_\xi$ using a large
rectangle. On the lower side, on the two vertical sides, and outside the
support of $g-1$, the symbol is $1$. On the upper side it is $g(x)$. The
positive orientation of the boundary traverses the upper side from
$x=+R$ to $x=-R$. Since $g$ has degree one when $x$ is traversed from
$-R$ to $+R$, the gluing along the boundary has degree $-1$. By
Definition~\ref{def:clase-bott}, the glued bundle is $\mathbf{E}_{-1}$ and
\[
 [\boldsymbol{\sigma}(T_b)]=[\mathbf{E}_{-1}]-[\underline{\mathbb C}]=b.
\]
\end{proof}

The preceding proof does not invoke a general Toeplitz index
formula: the only index calculated reduces to that of the unilateral shift.
With the conventions already specified, the argument agrees with the
normalization of
\cite[proof of Theorem~11.5, pp.~281--282]{BleeckerBoossIndex}.

\subsection{The oscillator and the Thom class}

On $L^2(\mathbb R)$ consider the closed operators
\begin{align}
 B^+u&=\frac{du}{dx}+xu,
 &B^-v&=-\frac{dv}{dx}+xv,
 \label{eq:osciladores-bott}
\end{align}
with common domain
\[
 \mathcal D
 :=\{u\in H^1(\mathbb R)\mid xu\in L^2(\mathbb R)\}.
\]
The space $C_c^\infty(\mathbb R)$ is dense in $\mathcal D$ for the norm
\[
 \|u\|_{\mathcal D}
 :=\|u\|_{L^2(\mathbb R)}+
 \left\|\frac{du}{dx}\right\|_{L^2(\mathbb R)}+
 \|xu\|_{L^2(\mathbb R)}.
\]
To justify density, fix $\chi\in C_c^\infty(\mathbb R)$,
$0\leq\chi\leq1$, equal to one on $[-1,1]$, and set
$\chi_R(x)=\chi(x/R)$. For $v\in\mathcal D$,
$\chi_Rv\to v$ in $L^2(\mathbb R)$ and
$x\chi_Rv\to xv$ in $L^2(\mathbb R)$ by dominated convergence.
Moreover,
\[
 (\chi_Rv)'-v'=(\chi_R-1)v'+\chi_R'v,
 \qquad
 \|\chi_R'v\|_{L^2(\mathbb R)}
 \leq R^{-1}\|\chi'\|_{L^\infty(\mathbb R)}\|v\|_{L^2(\mathbb R)},
\]
so $\chi_Rv\to v$ in $\mathcal D$.
For fixed $R$, convolution regularizations of $\chi_Rv$ converge
in $H^1(\mathbb R)$, and their supports remain in a fixed compact
interval. On that interval, multiplication by $x$ is bounded on
$L^2(\mathbb R)$, so they also converge in $\mathcal D$.
Choosing successively $R\to\infty$ and convolution radii producing
errors smaller than $1/R$ gives a sequence in $C_c^\infty(\mathbb R)$
converging in $\mathcal D$.

\begin{proposition}[Index of the Thom oscillator]
\label{prop:indice-oscilador-thom}
The operator $B^+:\mathcal D\longrightarrow L^2(\mathbb R)$ is Fredholm and
\[
 \ker B^+
 =\operatorname{span}\{e^{-\frac{x^2}{2}}\},
 \qquad
 \ker B^-=\{0\},
 \qquad
 \operatorname{ind}B^+=1.
\]
Its elliptic total symbol in the phase plane is $x+i\xi$ and represents
\[
 [\boldsymbol{\sigma}(B^+)]=-b=\lambda_{\mathbb C}.
\]
\end{proposition}

\begin{proof}
First we prove that
\begin{equation}
 \operatorname{Dom}\bigl((B^+)^*\bigr)=\mathcal D,
 \qquad
 (B^+)^*=B^-.
 \label{eq:dominio-adjunto-oscilador}
\end{equation}
Let $v\in\operatorname{Dom}((B^+)^*)$. There then exists
$f\in L^2(\mathbb R)$ such that
\[
 \int_{\mathbb R}
 \left(\frac{du}{dx}+xu\right)\overline v\,dx
 =\int_{\mathbb R}u\overline f\,dx,
 \qquad u\in C_c^\infty(\mathbb R).
\]
In the distributional sense,
\begin{equation}
 -\frac{dv}{dx}+xv=f.
 \label{eq:ecuacion-distribucional-adjunto-oscilador}
\end{equation}
Since $xv\in L^2_{\mathrm{loc}}(\mathbb R)$, this identity implies
$v\in H^1_{\mathrm{loc}}(\mathbb R)$.

Fix $\chi\in C_c^\infty(\mathbb R)$ with
$0\leq\chi\leq1$, $\chi=1$ on $[-1,1]$, and
$\operatorname{supp}\chi\subset[-2,2]$, and set
$\chi_R(x)=\chi\left(\frac{x}{R}\right)$ for $R\geq1$. The function
$w_R:=\chi_Rv$ belongs to $\mathcal D$ and
\begin{equation}
 B^-w_R=\chi_Rf-\chi_R'v.
 \label{eq:corte-adjunto-oscilador}
\end{equation}
For every compactly supported $w\in\mathcal D$, integration by parts
gives
\begin{equation}
 \|B^-w\|_{L^2(\mathbb R)}^2
 =\left\|\frac{dw}{dx}\right\|_{L^2(\mathbb R)}^2
  +\|xw\|_{L^2(\mathbb R)}^2
  +\|w\|_{L^2(\mathbb R)}^2.
 \label{eq:identidad-cortada-oscilador}
\end{equation}
Applied to $w_R$ and combined with
\eqref{eq:corte-adjunto-oscilador}, this equality yields the uniform bound
\begin{equation}
 \|x\chi_Rv\|_{L^2(\mathbb R)}^2
 \leq
 2\|f\|_{L^2(\mathbb R)}^2
 +\frac{2\|\chi'\|_{L^\infty(\mathbb R)}^2}{R^2}
  \|v\|_{L^2(\mathbb R)}^2.
 \label{eq:cota-cortes-dominio-adjunto}
\end{equation}
Fatou's lemma, applied as $R\to\infty$, shows that
$xv\in L^2(\mathbb R)$. The equation
\eqref{eq:ecuacion-distribucional-adjunto-oscilador} then gives
\[
 \frac{dv}{dx}=xv-f\in L^2(\mathbb R),
\]
so $v\in\mathcal D$. The reverse inclusion follows by integrating
by parts, now with both sides in $L^2(\mathbb R)$, and proves
\eqref{eq:dominio-adjunto-oscilador}.

The same argument shows that $(B^-)^*=B^+$ with domain $\mathcal D$.
Indeed, if $B^+u=g$ distributionally, then
\[
 \left\|\frac{dw}{dx}\right\|_{L^2(\mathbb R)}^2
 +\|xw\|_{L^2(\mathbb R)}^2
 =\|B^+w\|_{L^2(\mathbb R)}^2+
  \|w\|_{L^2(\mathbb R)}^2
\]
for $w=\chi_Ru$; the same passage to the limit gives $xu\in L^2(\mathbb R)$ and
$\frac{du}{dx}=g-xu\in L^2(\mathbb R)$. In particular, $B^+$ and $B^-$ are
closed operators and are adjoints of one another.

The homogeneous equations can be integrated directly:
\[
 B^+u=0
 \quad\Longleftrightarrow\quad
 u(x)=Ce^{-\frac{x^2}{2}},
\]
whereas
\[
 B^-v=0
 \quad\Longleftrightarrow\quad
 v(x)=Ce^{\frac{x^2}{2}}.
\]
Only the first solution belongs to $L^2(\mathbb R)$.

For $v\in C_c^\infty(\mathbb R)$, identity
\eqref{eq:identidad-cortada-oscilador} reads
\begin{align}
 \|B^-v\|_{L^2(\mathbb R)}^2
 &=\int_{\mathbb R}
 \left| -\frac{dv}{dx}+xv\right|^2\,dx\\
 &=\left\|\frac{dv}{dx}\right\|_{L^2(\mathbb R)}^2
   +\|xv\|_{L^2(\mathbb R)}^2
   +\|v\|_{L^2(\mathbb R)}^2.
 \label{eq:cota-oscilador-adjunto}
\end{align}
The approximation in the norm of $\mathcal D$ constructed before the
proposition preserves convergence of $v$, $v'$, and $xv$ simultaneously in
$L^2(\mathbb R)$. Since $B^-v=-v'+xv$, all terms in
\eqref{eq:cota-oscilador-adjunto} converge and the identity holds on
$\mathcal D$. In particular,
$B^-$ is bounded below by $1$ and has closed range. The closed
range theorem applied to the adjoint pair $B^+,B^-$ implies that
$\operatorname{Ran}B^+$ is closed. Moreover,
\[
 (\operatorname{Ran}B^+)^\perp=\ker B^-=\{0\},
\]
so $B^+$ is surjective. Its kernel has dimension one and its
cokernel is zero; this proves that its index is $1$.

With the Fourier normalization
$\widehat{\frac{du}{dx}}(\xi)=i\xi\widehat u(\xi)$, the total symbol is
$x+i\xi$. Its restriction to a large circle has degree one. By
the same gluing calculation as in the preceding proof, now with the usual
orientation of the circle, it represents $[E_1]-[\underline{\mathbb C}]=-b$.
Equality with $\lambda_{\mathbb C}$ is
\eqref{eq:thom-linea-menos-bott}.
\end{proof}

The oscillator does not belong to
$\operatorname{Ell}_c^0(\mathbb R)$: it is an unbounded operator. It does not have
compact resolvent; in fact,
$e^{\lambda x-\frac{x^2}{2}}\in\mathcal D$ and, for every
$\lambda\in\mathbb C$,
\[
 B^+e^{\lambda x-\frac{x^2}{2}}
 =\lambda e^{\lambda x-\frac{x^2}{2}}.
\]
The self-adjoint supersymmetric operator
\[
 \mathscr B=
 \begin{pmatrix}0&B^-\\ B^+&0\end{pmatrix}
\]
does have compact resolvent, since $\mathscr B^2$ is the sum of the two
scalar harmonic oscillators. The role of $B^+$ is to realize the
Thom class analytically, not to replace the model $T_b$. The comparison
\[
 \operatorname{ind}T_b=-1,
 \qquad
 \operatorname{ind}B^+=1
\]
confirms that $T_b$ represents $b$ and $B^+$ represents $-b$.

\section{Products of operators and multiplicativity of the index}
\label{sec:producto-operadores-indice}

Reducing dimension $n$ to the one-dimensional case requires checking
that external product of classes corresponds to a product of
operators and that their indices multiply. A subtlety cannot be
omitted: although operators acting separately in each variable
are well defined on the Hilbert tensor product, in general they are not
classical pseudodifferential operators on the product manifold. This phenomenon is
explained in \cite[Exercise~9.20, pp.~243--244]{BleeckerBoossIndex}. First we
calculate the index of the matrix on Hilbert spaces and then
regularize it using angular cutoffs.

Let
\[
 P:H_E^+\longrightarrow H_E^-,
 \qquad
 Q:H_F^+\longrightarrow H_F^-
\]
be bounded Fredholm operators between complex Hilbert spaces. On
the Hilbert tensor products define
\begin{equation}
 P\# Q
 :=
 \begin{pmatrix}
 P\otimes I&-I\otimes Q^*\\
 I\otimes Q&P^*\otimes I
 \end{pmatrix},
 \label{eq:producto-sharp-fredholm}
\end{equation}
with domain and codomain
\begin{align*}
 H^{\mathrm{ev}}
 &:=(H_E^+\widehat\otimes H_F^+)
   \oplus(H_E^-\widehat\otimes H_F^-),\\
 H^{\mathrm{odd}}
 &:=(H_E^-\widehat\otimes H_F^+)
   \oplus(H_E^+\widehat\otimes H_F^-).
\end{align*}

\begin{proposition}[Product of Fredholm operators]
\label{prop:multiplicatividad-indice-sharp}
The operator $P\#Q:H^{\mathrm{ev}}\longrightarrow H^{\mathrm{odd}}$ is
Fredholm and
\begin{equation}
 \operatorname{ind}(P\#Q)
 =\operatorname{ind}P\,\operatorname{ind}Q.
 \label{eq:indice-producto-sharp}
\end{equation}
\end{proposition}

\begin{proof}
Multiplying the matrices in \eqref{eq:producto-sharp-fredholm} and using
the fact that operators acting on different factors commute yields
\begin{align}
 (P\#Q)^*(P\#Q)
 &=
 \begin{pmatrix}
 P^*P\otimes I+I\otimes Q^*Q&0\\
 0&PP^*\otimes I+I\otimes QQ^*
 \end{pmatrix},
 \label{eq:cuadrado-sharp-par}\\
 (P\#Q)(P\#Q)^*
 &=
 \begin{pmatrix}
 PP^*\otimes I+I\otimes Q^*Q&0\\
 0&P^*P\otimes I+I\otimes QQ^*
 \end{pmatrix}.
 \label{eq:cuadrado-sharp-impar}
\end{align}
Each summand is nonnegative. If $\mathcal H_A,\mathcal H_B$ are Hilbert
spaces and $A\in\mathcal L(\mathcal H_A)$,
$B\in\mathcal L(\mathcal H_B)$ are nonnegative, then, on
$\mathcal H_A\widehat\otimes\mathcal H_B$,
\[
 \ker(A\otimes I+I\otimes B)=\ker A\widehat\otimes\ker B;
\]
indeed, for $u\in\mathcal H_A\widehat\otimes\mathcal H_B$,
\[
\begin{aligned}
 &\left\langle(A\otimes I+I\otimes B)u,u\right\rangle_{
 \mathcal H_A\widehat\otimes\mathcal H_B}\\
 &\qquad=
 \|(A^{\frac{1}{2}}\otimes I)u\|_{
 \mathcal H_A\widehat\otimes\mathcal H_B}^2
 +\|(I\otimes B^{\frac{1}{2}})u\|_{
 \mathcal H_A\widehat\otimes\mathcal H_B}^2.
\end{aligned}
\]
Since $P$ and $Q$ are Fredholm, their
restrictions to the orthogonal complements of their kernels are bounded
below. The sums in \eqref{eq:cuadrado-sharp-par} and
\eqref{eq:cuadrado-sharp-impar} are therefore bounded below on
the complements of their finite-dimensional kernels. This proves that
$P\#Q$ has closed range and finite-dimensional kernels.

The kernel formulas are
\begin{align*}
 \ker(P\#Q)
 &\cong
 (\ker P\widehat\otimes\ker Q)
 \oplus
 (\ker P^*\widehat\otimes\ker Q^*),\\
 \ker(P\#Q)^*
 &\cong
 (\ker P^*\widehat\otimes\ker Q)
 \oplus
 (\ker P\widehat\otimes\ker Q^*).
\end{align*}
If
\[
 a=\dim\ker P,
 \quad b_0=\dim\ker P^*,
 \quad c=\dim\ker Q,
 \quad d=\dim\ker Q^*,
\]
then
\[
 \operatorname{ind}(P\#Q)
 =(ac+b_0d)-(b_0c+ad)
 =(a-b_0)(c-d),
\]
which proves \eqref{eq:indice-producto-sharp}.
\end{proof}

The preceding matrix cannot simply be declared classical pseudodifferential on a
product. To obtain a classical matrix while preserving its index, we first raise
the order. Let $(X,g_X)$ and $(Y,g_Y)$ be closed Riemannian manifolds and let
\[
 P\in\Psi^0_{\mathrm{cl}}(X;\mathbb C^q),
 \qquad
 Q\in\Psi^0_{\mathrm{cl}}(Y;\mathbb C^r)
\]
be elliptic. Choose nonnegative scalar Laplacians $\Delta_X$ and
$\Delta_Y$ and
set
\[
 \Lambda_{X,q}=(I+\Delta_X)^{\frac{1}{2}}\otimes I_{\mathbb C^q},
 \qquad
 \Lambda_{Y,r}=(I+\Delta_Y)^{\frac{1}{2}}\otimes I_{\mathbb C^r}.
\]
By Theorem~\ref{pseudo:teo-potencias-complejas-clasicas}, these are positive classical
operators of order one and define isomorphisms
\[
 \Lambda_{X,q}:H^1(X;\mathbb C^q)\longrightarrow L^2(X;\mathbb C^q),
 \qquad
 \Lambda_{Y,r}:H^1(Y;\mathbb C^r)\longrightarrow L^2(Y;\mathbb C^r).
\]
The operators
\begin{equation}
 A:=P\Lambda_{X,q}\in\Psi^1_{\mathrm{cl}}(X;\mathbb C^q),
 \qquad
 B:=Q\Lambda_{Y,r}\in\Psi^1_{\mathrm{cl}}(Y;\mathbb C^r)
 \label{eq:producto-operadores-orden-uno}
\end{equation}
are elliptic of order one and
\begin{equation}
 \operatorname{ind}A=\operatorname{ind}P,
 \qquad
 \operatorname{ind}B=\operatorname{ind}Q.
 \label{eq:indices-tras-elevar-orden}
\end{equation}

Define
\begin{equation}
 \mathcal D\colon
 C^\infty(X\times Y;\mathbb C^{qr})^{\oplus2}
 \longrightarrow C^\infty(X\times Y;\mathbb C^{qr})^{\oplus2},
 \qquad
 \mathcal D:=A\#B
 =\begin{pmatrix}
 A\otimes I_{\mathbb C^r}&-I_{\mathbb C^q}\otimes B^*\\
 I_{\mathbb C^q}\otimes B&A^*\otimes I_{\mathbb C^r}
 \end{pmatrix}.
 \label{eq:producto-no-regularizado-orden-uno}
\end{equation}
Its closure, as an operator on
$L^2(X\times Y;\mathbb C^{qr})^{\oplus2}$, has domain
$H^1(X\times Y;\mathbb C^{qr})^{\oplus2}$. Indeed, if
$\mathcal L_X:=\Lambda_{X,q}\otimes I_{\mathbb C^r}$ and
$\mathcal L_Y:=I_{\mathbb C^q}\otimes\Lambda_{Y,r}$, then
\begin{equation}
 \mathcal L
 :=(\mathcal L_X^2+\mathcal L_Y^2)^{\frac{1}{2}}
 \colon H^1(X\times Y;\mathbb C^{qr})
 \longrightarrow L^2(X\times Y;\mathbb C^{qr})
 \label{eq:operador-sobolev-producto}
\end{equation}
is an isomorphism, and the same holds for $\mathcal L^{\oplus2}$. The
elliptic estimates for $A$ and $B$ show that the graph norm of
$\mathcal D$ is equivalent to that of
$H^1(X\times Y;\mathbb C^{qr})^{\oplus2}$.

The calculation in Proposition~\ref{prop:multiplicatividad-indice-sharp}
remains valid for these closed operators. The squares
$A^*A$, $AA^*$, $B^*B$, and $BB^*$ now have compact resolvent. The spectral
theorem applied to
\begin{align}
 \mathcal D^*\mathcal D
 &=\begin{pmatrix}
 A^*A\otimes I+I\otimes B^*B&0\\
 0&AA^*\otimes I+I\otimes BB^*
 \end{pmatrix},
 \label{eq:cuadrado-producto-orden-uno-par}\\
 \mathcal D\mathcal D^*
 &=\begin{pmatrix}
 AA^*\otimes I+I\otimes B^*B&0\\
 0&A^*A\otimes I+I\otimes BB^*
 \end{pmatrix}
 \label{eq:cuadrado-producto-orden-uno-impar}
\end{align}
gives the same four kernel formulas as in the preceding proof. In
particular,
\begin{equation}
 \operatorname{ind}\mathcal D
 =\operatorname{ind}A\,\operatorname{ind}B
 =\operatorname{ind}P\,\operatorname{ind}Q.
 \label{eq:indice-producto-no-regularizado}
\end{equation}

The approximation we need is constructed within the ordinary
pseudodifferential calculus on the product. In this way we do not assume
that an operator defined by spectral calculus in the two variables is
automatically an isotropic pseudodifferential operator.

\begin{lemma}[Angular microlocal cutoffs]
\label{lem:cortes-microlocales-angulares}
Fix finite atlases $\{U_i\}_{i\in I}$ of $X$ and $\{V_j\}_{j\in J}$ of $Y$, together with
real functions
\[
 \alpha_i\in C_c^\infty(U_i),
 \qquad
 \beta_j\in C_c^\infty(V_j),
 \qquad
 \sum_{i\in I}\alpha_i^2=1,
 \qquad
 \sum_{j\in J}\beta_j^2=1.
\]
There is a choice, depending only on these data, of scalar operators
\[
 Q_R^X,Q_R^Y\in\Psi^0_{\mathrm{cl}}(X\times Y),
 \qquad R>2,
\]
with the following property. If
$T\in\Psi^1_{\mathrm{cl}}(X;\mathbb C^d)$ and
$S\in\Psi^1_{\mathrm{cl}}(Y;\mathbb C^e)$, extend
$Q_R^X,Q_R^Y$ diagonally to $\mathbb C^{de}$ and set
\begin{equation}
 \mathfrak R_R^X(T):=(T\otimes I_{\mathbb C^e})Q_R^X,
 \qquad
 \mathfrak R_R^Y(S):=(I_{\mathbb C^d}\otimes S)Q_R^Y
 \label{eq:def-regularizaciones-microlocales}
\end{equation}
These operators belong to
$\Psi^1_{\mathrm{cl}}(X\times Y;\mathbb C^{de})$ and satisfy
\begin{align}
 \bigl\|\mathfrak R_R^X(T)-T\otimes I_{\mathbb C^e}\bigr\|_{
 \mathcal L(H^1(X\times Y;\mathbb C^{de}),
 L^2(X\times Y;\mathbb C^{de}))}
 &\leq\frac{C_T}{R},
 \label{eq:cota-regularizacion-microlocal-X}\\
 \bigl\|\mathfrak R_R^Y(S)-I_{\mathbb C^d}\otimes S\bigr\|_{
 \mathcal L(H^1(X\times Y;\mathbb C^{de}),
 L^2(X\times Y;\mathbb C^{de}))}
 &\leq\frac{C_S}{R}.
 \label{eq:cota-regularizacion-microlocal-Y}
\end{align}

Let $|\xi|_i$ be the Euclidean norm of the components of
$\xi\in T_x^*X$ in the chart $U_i$, and define $|\eta|_j$ similarly.
In each summand of the following formulas, these norms are used only on
$\operatorname{supp}\alpha_i\times\operatorname{supp}\beta_j$. The principal
symbols of $T$ and $S$ are denoted, respectively, by
\[
 t_1(x,\xi):=\boldsymbol{\sigma}_1^\Psi(T)(x,\xi),
 \qquad
 s_1(y,\eta):=\boldsymbol{\sigma}_1^\Psi(S)(y,\eta),
\]
and are understood to be lifted to $T^*(X\times Y)$ in the region where their
own covariable does not vanish. The principal symbols of the two
cutoff operators are
\begin{align}
 \Phi_R(x,y,\xi,\eta)
 &:=\sum_{(i,j)\in I\times J}\alpha_i(x)^2\beta_j(y)^2
 \chi\left(\frac{|\eta|_j}{R|\xi|_i}\right),
 \label{eq:simbolo-principal-regularizacion-X}\\
 \Psi_R(x,y,\xi,\eta)
 &:=\sum_{(i,j)\in I\times J}\alpha_i(x)^2\beta_j(y)^2
 \chi\left(\frac{|\xi|_i}{R|\eta|_j}\right).
 \label{eq:simbolo-principal-regularizacion-Y}
\end{align}
Here $\chi\in C^\infty([0,\infty),[0,1])$ equals $1$ on $[0,1]$, is
positive on $[0,2)$, and equals $0$ on $[2,\infty)$. On the axis $\xi=0$, the
summand of \eqref{eq:simbolo-principal-regularizacion-X} is defined to be
$0$, and on the axis $\eta=0$ to be $1$; for
\eqref{eq:simbolo-principal-regularizacion-Y} the opposite convention is used.
In particular,
\begin{equation}
 \boldsymbol{\sigma}_1^\Psi\bigl(\mathfrak R_R^X(T)\bigr)=\Phi_Rt_1,
 \qquad
 \boldsymbol{\sigma}_1^\Psi\bigl(\mathfrak R_R^Y(S)\bigr)=\Psi_Rs_1.
 \label{eq:simbolos-principales-regularizaciones}
\end{equation}
The products in \eqref{eq:simbolos-principales-regularizaciones} extend
smoothly by zero to the axis of the corresponding covariable: the
factor $\Phi_R$ is identically zero on a conic neighborhood of $\xi=0$ when
$\eta\neq0$, and $\Psi_R$ has the symmetric property.
\end{lemma}

\begin{proof}
Set $\lambda_{ij}(x,y):=\alpha_i(x)\beta_j(y)$. Identify
$U_i\times V_j$ with an open subset of $\mathbb R^{m+k}$ and write
\[
 \langle\xi\rangle_i=(1+|\xi|_i^2)^{\frac{1}{2}},
 \qquad
 \langle\eta\rangle_j=(1+|\eta|_j^2)^{\frac{1}{2}}.
\]
Let $M_{ij,R}^X$ and $M_{ij,R}^Y$ be the Fourier multipliers in this
product chart with full symbols
\begin{align}
 \vartheta_{ij,R}^X(\xi,\eta)
 &:=\chi\left(
 \frac{\langle\eta\rangle_j}{R\langle\xi\rangle_i}
 \right),
 \label{eq:corte-angular-completo-carta-X}\\
 \vartheta_{ij,R}^Y(\xi,\eta)
 &:=\chi\left(
 \frac{\langle\xi\rangle_i}{R\langle\eta\rangle_j}
 \right).
 \label{eq:corte-angular-completo-carta-Y}
\end{align}
Multiplication by $\lambda_{ij}$ before and after the multiplier makes
its extension by zero outside the chart well defined. Take
\begin{equation}
 Q_R^X:=\sum_{(i,j)\in I\times J}\lambda_{ij}M_{ij,R}^X\lambda_{ij},
 \qquad
 Q_R^Y:=\sum_{(i,j)\in I\times J}\lambda_{ij}M_{ij,R}^Y\lambda_{ij}.
 \label{eq:operadores-corte-angular-atlas}
\end{equation}
If all multipliers in \eqref{eq:operadores-corte-angular-atlas} are replaced
by the identity, the result is exactly the identity,
since
\begin{equation}
 \sum_{(i,j)\in I\times J}\lambda_{ij}^2
 =\left(\sum_{i\in I}\alpha_i^2\right)
  \left(\sum_{j\in J}\beta_j^2\right)=1.
 \label{eq:particion-cuadratica-producto}
\end{equation}

For fixed $R$, the symbols of
\eqref{eq:corte-angular-completo-carta-X} and
\eqref{eq:corte-angular-completo-carta-Y} are classical of order zero. We
verify this for the first. If $\xi\neq0$ and $\eta\neq0$, then
\begin{equation}
 \frac{\langle t\eta\rangle_j}{R\langle t\xi\rangle_i}
 \sim
 \frac{|\eta|_j}{R|\xi|_i}
 \left(1+\sum_{\ell=1}^{\infty}c_{\ell,ij}(\xi,\eta)t^{-2\ell}\right),
 \qquad t\longrightarrow\infty.
 \label{eq:expansion-radial-corte-angular}
\end{equation}
Taylor's formula for $\chi$ gives a homogeneous expansion of degrees
$0,-2,-4,\ldots$. Its coefficients extend smoothly to the axes:
near $\{\xi=0,\ \eta\neq0\}$ the cutoff is identically zero, and near
$\{\eta=0,\ \xi\neq0\}$ it is identically one. The homogeneous term of degree
zero is
\[
 \chi\left(\frac{|\eta|_j}{R|\xi|_i}\right).
\]
The same proof, with variables interchanged, applies to the second
symbol. Thus $Q_R^X$ and $Q_R^Y$ are classical operators of order zero,
with principal symbols $\Phi_R$ and $\Psi_R$, respectively.

We must check that the compositions in
\eqref{eq:def-regularizaciones-microlocales} belong to the isotropic
calculus; this does not follow merely from $Q_R^X$ and $Q_R^Y$ being
pseudodifferential operators. Consider $Q_R^X$. In a chart, the symbol of every
summand of \eqref{eq:operadores-corte-angular-atlas}, modulo
$S^{-\infty}$, has an expansion whose terms contain a derivative of
$\vartheta_{ij,R}^X$. All these terms are supported in
\begin{equation}
 \langle\eta\rangle_j\leq2R\langle\xi\rangle_i,
 \qquad
 \langle(\xi,\eta)\rangle
 \leq C_R\langle\xi\rangle_i.
 \label{eq:comparacion-pesos-cono-X}
\end{equation}
Asymptotic summation can be performed after multiplication by a cutoff
equal to $1$ on this cone; it therefore retains support in a slightly
larger cone, up to a remainder in $S^{-\infty}$.

Use a finite partition to decompose the kernel of $T$ into a part
near the diagonal and a smooth part. For the first part, the composition formula
contains terms
\begin{equation}
 \frac{1}{\gamma!}
 \partial_\xi^\gamma t(x,\xi)
 D_x^\gamma q_{ij,R}(x,y,\xi,\eta),
 \label{eq:expansion-simbolica-corte-angular}
\end{equation}
where $q_{ij,R}$ is a full symbol of the corresponding summand of
$Q_R^X$. On the cone \eqref{eq:comparacion-pesos-cono-X}, the estimates for
$t$ with respect to $\xi$ become isotropic estimates of order
one with respect to $(\xi,\eta)$. The remainder formula, integrated by parts in
the chart variables, gives the same estimates. The smooth part of the
kernel of $T$ has rapidly decreasing transform in $\xi$; since
$|\eta|_j\leq C_R\langle\xi\rangle_i$ on the support of the cutoff, its
composition is rapidly decreasing in $(\xi,\eta)$. This proves that
$(T\otimes I_{\mathbb C^e})Q_R^X$ is classical of order one. Moreover, the first term in
\eqref{eq:expansion-simbolica-corte-angular} proves the first equality in
\eqref{eq:simbolos-principales-regularizaciones}. The proof for
$(I_{\mathbb C^d}\otimes S)Q_R^Y$ is symmetric.

The norm estimate remains. From
\eqref{eq:particion-cuadratica-producto} we obtain the exact identity
\begin{equation}
 Q_R^X-I
 =\sum_{(i,j)\in I\times J}\lambda_{ij}(M_{ij,R}^X-I)\lambda_{ij}.
 \label{eq:diferencia-corte-angular-identidad}
\end{equation}
Let $u\in H^1(X\times Y;\mathbb C^{de})$ and let $v$ be the coordinate
expression, with values in $\mathbb C^{de}$, of $\lambda_{ij}u$. The factor
$1-\vartheta_{ij,R}^X$ can be nonzero only when
$\langle\eta\rangle_j>R\langle\xi\rangle_i$. By Plancherel,
\begin{align}
 &\bigl\|(M_{ij,R}^X-I)v\bigr\|_{H^1(\mathbb R^m_x;
 L^2(\mathbb R^k_y;\mathbb C^{de}))}^2\notag\\
 &\quad=
 \int_{\mathbb R^k}\int_{\mathbb R^m}
 \langle\xi\rangle_i^2
 |1-\vartheta_{ij,R}^X(\xi,\eta)|^2
 |\widehat v(\xi,\eta)|_{\mathbb C^{de}}^2\,d\xi\,d\eta\notag\\
 &\quad\leq\frac{1}{R^2}
 \int_{\mathbb R^k}\int_{\mathbb R^m}
 \langle\eta\rangle_j^2
 |\widehat v(\xi,\eta)|_{\mathbb C^{de}}^2\,d\xi\,d\eta
 =\frac{1}{R^2}
 \|v\|_{L^2(\mathbb R^m_x;
 H^1(\mathbb R^k_y;\mathbb C^{de}))}^2.
 \label{eq:plancherel-corte-angular}
\end{align}
Suppressing the chart identifications from the notation, we have
\begin{equation}
 (T\otimes I_{\mathbb C^e})(Q_R^X-I)u
 =\sum_{(i,j)\in I\times J}(T\alpha_i\otimes\beta_j I_{\mathbb C^e})
 (M_{ij,R}^X-I)(\lambda_{ij}u).
 \label{eq:diferencia-regularizada-suma-local}
\end{equation}
The operator $T\alpha_i:H^1(X;\mathbb C^d)\to L^2(X;\mathbb C^d)$ is
bounded. Multiplications
by $\alpha_i$ and $\beta_j$, changes of coordinates, and
their inverses are bounded on the Sobolev spaces appearing in
\eqref{eq:plancherel-corte-angular}. Applying $T\alpha_i$ for every value of
$y$ and using Fubini yields
\[
 \|(T\alpha_i\otimes\beta_j I_{\mathbb C^e})w\|_{
 L^2(X\times Y;\mathbb C^{de})}
 \leq C_{ij}\|w\|_{H^1(\mathbb R^m_x;
 L^2(\mathbb R^k_y;\mathbb C^{de}))}.
\]
Since the covers are finite,
\eqref{eq:diferencia-corte-angular-identidad} and
\eqref{eq:plancherel-corte-angular} give
\eqref{eq:cota-regularizacion-microlocal-X}, with a constant independent
of $R$. Interchanging $X$ and $Y$ gives
\eqref{eq:cota-regularizacion-microlocal-Y}.
\end{proof}

The two principal symbols constructed with the same atlas do not vanish
simultaneously. Indeed, every summand of $\Phi_R$ and $\Psi_R$ is nonnegative,
and the weights $\alpha_i^2\beta_j^2$ sum to $1$. First suppose
that $\xi\neq0$ and $\eta\neq0$. If $\Phi_R=0$, then for every active pair we have
\[
 \frac{|\eta|_j}{R|\xi|_i}\geq2.
\]
Therefore,
\[
 \frac{|\xi|_i}{R|\eta|_j}
 \leq\frac{1}{2R^2}<1
\]
for every active pair. All cutoffs occurring in $\Psi_R$ equal
$1$, whence $\Psi_R=1$. The same reasoning with the variables
interchanged shows that $\Psi_R=0$ implies $\Phi_R=1$. On the axes, the
conclusion follows directly from the conventions in
\eqref{eq:simbolo-principal-regularizacion-X} and
\eqref{eq:simbolo-principal-regularizacion-Y}.

Apply the lemma to $A$, $A^*$, $B$, and $B^*$ and define
\begin{equation}
 \mathcal D_R
 :=\begin{pmatrix}
 \mathfrak R_R^X(A)&-\mathfrak R_R^Y(B^*)\\
 \mathfrak R_R^Y(B)&\mathfrak R_R^X(A^*)
 \end{pmatrix}.
 \label{eq:producto-regularizado-angular}
\end{equation}
This is an operator
\[
 \mathcal D_R\in
 \Psi^1_{\mathrm{cl}}(X\times Y;\mathbb C^{2qr}).
\]
The four estimates in the lemma give
\begin{equation}
 \|\mathcal D_R-\mathcal D\|_{
 \mathcal L(
 H^1(X\times Y;\mathbb C^{qr})^{\oplus2},
 L^2(X\times Y;\mathbb C^{qr})^{\oplus2})}
 \leq\frac{C}{R}.
 \label{eq:convergencia-norma-producto-regularizado}
\end{equation}
Since $\mathcal D$ is Fredholm, openness of the set of Fredholm
operators implies that, for $R$ sufficiently large,
\begin{equation}
 \operatorname{ind}\mathcal D_R
 =\operatorname{ind}\mathcal D
 =\operatorname{ind}P\,\operatorname{ind}Q.
 \label{eq:indice-producto-regularizado}
\end{equation}

It remains to check that the symbol of $\mathcal D_R$ is the product complex.
Let $a_1(x,\xi)$ and $b_1(y,\eta)$ be the principal symbols of $A$ and $B$.
The full principal symbol is
\begin{equation}
 c_R(x,y,\xi,\eta)
 =\begin{pmatrix}
 \Phi_R a_1\otimes I&-\Psi_R I\otimes b_1^*\\
 \Psi_R I\otimes b_1&\Phi_R a_1^*\otimes I
 \end{pmatrix}.
 \label{eq:simbolo-producto-regularizado}
\end{equation}
Matrix multiplication gives
\begin{align}
 c_R^*c_R
 =\begin{pmatrix}
 \Phi_R^2a_1^*a_1\otimes I+
  \Psi_R^2I\otimes b_1^*b_1&0\\
 0&\Psi_R^2I\otimes b_1b_1^*+
  \Phi_R^2a_1a_1^*\otimes I
 \end{pmatrix}.
 \label{eq:positividad-simbolo-producto-regularizado}
\end{align}
Let $(\xi,\eta)\neq(0,0)$. If $\xi\neq0$ and $\Phi_R>0$, ellipticity of
$a_1$ makes the contribution containing $a_1^*a_1$ in the
first diagonal entry and that containing $a_1a_1^*$ in the second positive definite. If
$\xi=0$ or $\Phi_R=0$, then $\eta\neq0$ and $\Psi_R=1$ by the preceding
property; ellipticity of $b_1$ makes the remaining contributions
positive definite. Thus both diagonal entries of
\eqref{eq:positividad-simbolo-producto-regularizado} are positive definite
on the entire punctured cotangent bundle. Compactness of the cosphere bundle of
$X\times Y$ gives a uniform lower bound, so $c_R$ is elliptic.

To identify the class, a continuous homotopy of exact complexes
over the cosphere bundle suffices, as this is one of the relations defining $K_c^0$;
we do not assert that its intermediate terms are classical symbols.
Extend $a_1$ and $b_1$ continuously by zero on the axes of their
covariables. For $t<1$ these extensions may be merely continuous on
the axes, which suffices for this topological homotopy. Set
\[
 \Phi_{R,t}:=(1-t)+t\Phi_R,
 \qquad
 \Psi_{R,t}:=(1-t)+t\Psi_R,
 \qquad 0\leq t\leq1,
\]
and replace $\Phi_R,\Psi_R$ by
$\Phi_{R,t},\Psi_{R,t}$ in \eqref{eq:simbolo-producto-regularizado}.
For $t<1$ both coefficients are positive. For $t=1$ they do not vanish
simultaneously. The calculation in
\eqref{eq:positividad-simbolo-producto-regularizado} proves that the entire
continuous family of complexes is invertible away from the zero section. At
$t=0$ we obtain
\begin{equation}
 a_1\#b_1
 =\begin{pmatrix}
 a_1\otimes I&-I\otimes b_1^*\\
 I\otimes b_1&a_1^*\otimes I
 \end{pmatrix},
 \label{eq:complejo-producto-sin-cortes}
\end{equation}
which is the tensor product complex of the two complexes of length one.
This matrix expresses the sign rule for the product of complexes described in
\cite[Remark~11.3, pp.~278--279]{BleeckerBoossIndex}.
Proposition~\ref{prop:producto-externo-k} and homotopy invariance
therefore give
\begin{equation}
 [c_R]=[\boldsymbol{\sigma}(P)]\boxtimes[\boldsymbol{\sigma}(Q)].
 \label{eq:clase-producto-regularizado}
\end{equation}
Multiplication by the positive symbols of $\Lambda_X$ and
$\Lambda_Y$ does not change the classes of $P$ and $Q$. This argument is the
Euclidean version of the cutoffs $c_{r_0}$ used in the proof of the
multiplicative property in
\cite[pp.~294--301]{BleeckerBoossIndex}; raising the order to one beforehand
also makes possible the norm comparison
\eqref{eq:convergencia-norma-producto-regularizado}.

Finally,
\begin{equation}
 C_R:=\mathcal D_R(\mathcal L^{-1})^{\oplus2}
 \in\Psi^0_{\mathrm{cl}}(X\times Y;\mathbb C^{2qr})
 \label{eq:producto-regularizado-orden-cero}
\end{equation}
has the same index as $\mathcal D_R$. Its symbol is obtained by dividing
$c_R$ by the positive function
$\sqrt{|\xi|_{\mathbf{g}_X}^2+|\eta|_{\mathbf{g}_Y}^2}$; this division also leaves the class
\eqref{eq:clase-producto-regularizado} unchanged.

Apply the construction to
$P\in\operatorname{Ell}_c^0(\mathbb R^m)$ and
$Q\in\operatorname{Ell}_c^0(\mathbb R^k)$. Extend them by the identity
to $X=S^m$ and $Y=S^k$. The class on the right-hand side of
\eqref{eq:clase-producto-regularizado} is compactly supported inside
\[
 T^*(\mathbb R^m\times\mathbb R^k)
 =T^*\mathbb R^{m+k}.
\]
The relative realization of
Proposition~\ref{prop:indice-analitico-euclidiano-k} allows us to stabilize the
complex, deform it while keeping it fixed on a compact set to the identity near
the divisor
\[
 (\{\infty\}\times S^k)\cup(S^m\times\{\infty\})
\]
and quantize the homotopy. If $\rho$ equals $1$ on the projection of the
symbol support, the final replacement
\[
 \widetilde C\longmapsto
 I+\rho(\widetilde C-I)\rho
\]
produces an operator in
$\operatorname{Ell}_c^0(\mathbb R^{m+k})$ with the same symbol class.
The intermediate operators are elliptic on the closed manifold
$S^m\times S^k$; by
Theorem~\ref{teo:indice-analitico-sobre-k}, their index depends only on this
class. Combined with \eqref{eq:indice-producto-regularizado}, this
proves
\begin{equation}
 I_{m+k}(a\boxtimes c)=I_m(a)I_k(c),
 \qquad
 a\in K_c^0(T^*\mathbb R^m),
 \quad c\in K_c^0(T^*\mathbb R^k).
 \label{eq:multiplicatividad-indice-euclidiano-k}
\end{equation}

\section{The index theorem in Euclidean space}
\label{sec:teorema-indice-euclidiano}

Define the Euclidean topological index by
\begin{equation}
 \operatorname{ind}_{\mathrm t}^{\mathbb R^n}
 :=(-1)^n\alpha_n:
 K_c^0(T^*\mathbb R^n)=K_c^0(\mathbb R^{2n})
 \longrightarrow\mathbb Z.
 \label{eq:def-indice-topologico-euclidiano}
\end{equation}
Equality \eqref{eq:inversa-pushforward-punto} shows that it can also
be written as
\[
 \operatorname{ind}_{\mathrm t}^{\mathbb R^n}=i_!^{-1},
 \qquad i:\{0\}\hookrightarrow\mathbb R^n.
\]

\begin{theorem}[Euclidean index]
\label{teo:indice-euclidiano-k}
Let $P\in\operatorname{Ell}_c^0(\mathbb R^n)$. Then
\begin{equation}
 \operatorname{ind}_{\mathrm a}(P)
 =\operatorname{ind}_{\mathrm t}^{\mathbb R^n}[\boldsymbol{\sigma}(P)]
 =(-1)^n\alpha_n([\boldsymbol{\sigma}(P)]).
 \label{eq:formula-indice-euclidiano}
\end{equation}
\end{theorem}

\begin{proof}
By Corollary~\ref{cor:iteracion-bott}, every class in
$K_c^0(T^*\mathbb R^n)$ has a unique expression
\[
 [\boldsymbol{\sigma}(P)]=\ell\,b^{\boxtimes n},
 \qquad \ell\in\mathbb Z,
\]
and
\[
 \ell=\alpha_n([\boldsymbol{\sigma}(P)]).
\]
Lemma~\ref{lem:modelo-bott-indice} and multiplicativity
\eqref{eq:multiplicatividad-indice-euclidiano-k} imply
\[
 I_n(b^{\boxtimes n})
 =I_1(b)^n=(-1)^n.
\]
Since $I_n$ is a homomorphism,
\[
 \operatorname{ind}_{\mathrm a}(P)
 =I_n([\boldsymbol{\sigma}(P)])
 =\ell I_n(b^{\boxtimes n})
 =(-1)^n\alpha_n([\boldsymbol{\sigma}(P)]).
\]
By \eqref{eq:def-indice-topologico-euclidiano}, the last expression is the
topological index.
\end{proof}

The proof uses Bott periodicity to identify the generator and
calculates its index using the unilateral shift and angular
regularization; it therefore does not presuppose the index theorem on closed manifolds.

\subsection{An extension: the families index}

Let $\Lambda$ be a compact Hausdorff space and let
$P=(P_\lambda)_{\lambda\in\Lambda}$ be a family of Euclidean elliptic
operators of order zero, continuous in the topology of classical symbols and
with a common compact support. In particular, both the operators and
their principal symbols form continuous families. Their analytic index is
a class in $K^0(\Lambda)$. Recall its
construction. Choose a finite-dimensional subspace
$V\subset L^2(\mathbb R^n;\mathbb C^q)$ such that
\[
 \operatorname{Ran}P_\lambda+V
 =L^2(\mathbb R^n;\mathbb C^q)
 \qquad\text{for every }\lambda\in\Lambda.
\]
Existence of a single $V$ for all parameters is obtained first
locally from openness
of surjectivity and then by taking the sum of the spaces
corresponding to a finite subcover of $\Lambda$. The kernels of
the surjective maps
\[
 \mathcal P_\lambda:
 L^2(\mathbb R^n;\mathbb C^q)\oplus V
 \longrightarrow L^2(\mathbb R^n;\mathbb C^q),
 \qquad
 \mathcal P_\lambda(u,v)=P_\lambda u+v,
\]
form a finite-dimensional vector bundle
$\mathbf{K}:=\ker\mathcal P\longrightarrow\Lambda$. Indeed,
fix $\lambda_0$ and let $R_0$ be a right inverse of
$\mathcal P_{\lambda_0}$. For $\lambda$ near $\lambda_0$, the operator
$\mathcal P_\lambda R_0$ is invertible by the Neumann series, and
\[
 R_\lambda
 :=R_0(\mathcal P_\lambda R_0)^{-1}
\]
is a right inverse depending continuously on $\lambda$. The
projections
\[
 \Pi_\lambda:=I-R_\lambda\mathcal P_\lambda
\]
have image $\ker\mathcal P_\lambda$, depend continuously on the parameter,
and provide local trivializations of these kernels. Define
\begin{equation}
 \operatorname{Ind}_{\mathrm a}(P)
 :=[\mathbf{K}]-[\underline V]
 \in K^0(\Lambda).
 \label{eq:indice-analitico-familia-euclidiana}
\end{equation}
The class is independent of $V$. Indeed, if $V\subseteq W$ and both
spaces make the stabilized
maps surjective, set
$\mathbf{K}_V:=\ker\mathcal P_V$ and
$\mathbf{K}_W:=\ker\mathcal P_W$. There is an exact sequence of bundles
\[
 0\longrightarrow\mathbf{K}_V
 \longrightarrow\mathbf{K}_W
 \longrightarrow\underline{W/V}
 \longrightarrow0.
\]
The last arrow sends $(u,w)$ to the class of $w$ modulo $V$. It is
surjective because
$\operatorname{Ran}P_\lambda+V=L^2(\mathbb R^n;\mathbb C^q)$, and its kernel
is $\mathbf{K}_V$. Thus
\[
 [\mathbf{K}_W]-[\underline W]
 =[\mathbf{K}_V]+[\underline{W/V}]
  -[\underline V]-[\underline{W/V}]
 =[\mathbf{K}_V]-[\underline V].
\]
Two arbitrary choices are compared by passing to their sum.

\begin{proposition}[Realization and descent for families]
\label{prop:realizacion-descenso-familias-euclidianas}
Let $\Lambda$ be a compact Hausdorff space. Every class
\[
 \kappa\in K_c^0(T^*\mathbb R^n\times\Lambda)
\]
is the symbol class of a family of operators in
$\operatorname{Ell}_c^0(\mathbb R^n)$ continuous in the symbol topology,
after a fixed stabilization. The index bundle of a family depends
only on its symbol class. Consequently, there exists a well-defined
homomorphism
\begin{equation}
 I_{n,\Lambda}:
 K_c^0(\mathbb R^{2n}\times\Lambda)
 \longrightarrow K^0(\Lambda)
 \label{eq:indice-familiar-sobre-k}
\end{equation}
sending a class realized by $P$ to
$\operatorname{Ind}_{\mathrm a}(P)$.
\end{proposition}

\begin{proof}
Extend $\kappa$ by zero to $T^*S^n\times\Lambda$. The projection of
its compact support onto $T^*\mathbb R^n$ is compact. We may choose
a common neighborhood $U$ of $\infty$ whose closure does not meet the projection
of this support onto $S^n$. Under the radial identification of $T^*S^n$
with the interior of $B^*S^n$, the support lies in the complement of the
closed set
\begin{equation}
 \bigl(S^*S^n\times\Lambda\bigr)
 \cup\bigl(B^*\overline U\times\Lambda\bigr).
 \label{eq:subespacio-relativo-familias}
\end{equation}
The description by supports and excision in
Theorem~\ref{teo:exactitud-excision-k-apendice} allow the
class to be represented in the relative group of this pair. Since the total space and the
closed subspace are compact Hausdorff,
Theorem~\ref{teo:modelo-triples-k-relativa} provides bundles and an
isomorphism over this closed subset, elementary on the second set.
No CW structure on $\Lambda$ is needed. The relative isomorphism
extends to a neighborhood of the closed subset by the extension argument
in the proof of that theorem. In particular, the relative
trivialization is available near $B^*\overline U\times\Lambda$.
Radial contraction of $B^*S^n$ reduces
the bundles, after stabilization, to pullbacks of bundles over
$S^n\times\Lambda$. On restriction to
$\mathbb R^n\times\Lambda$, these bundles are pullbacks of bundles $G$ and $H$
over $\Lambda$. The elementary trivialization near $\infty$ identifies
$G$ with $H$. Choose a complement $G'$ such that
$G\oplus G'\cong\underline{\mathbb C}^{\,N}$ and add the elementary triple
defined by $G'$. Thus, over $\mathbb R^n\times\Lambda$ the operators
act on a trivial bundle of fixed rank and are the identity on the added
summands.

Smoothing must be performed uniformly. Take a finite atlas of
$S^n$ and a partition of unity independent of $\lambda$, and
regularize by convolution only in the variables of $S^n$ and their
covariables. Shrink $U$ once and use a cutoff function
vanishing near its closure to keep the relative trivialization
intact; continue to denote the smaller neighborhood by $U$. The original isomorphism
on the compact set $S^*S^n\times\Lambda$ satisfies
\[
 \delta:=\inf_{(x,\xi,\lambda)\in S^*S^n\times\Lambda}
 \bigl\|\theta(x,\xi,\lambda)^{-1}\bigr\|_{
 \operatorname{op}(\mathbb C^N)}^{-1}>0.
\]
Denote the approximation by $\widetilde\theta$ and choose it so that
\[
 \|\widetilde\theta-\theta\|_{
 C^0(S^*S^n\times\Lambda,\operatorname{End}(\mathbb C^N))}
 <\frac{\delta}{2}.
\]
The segment between the two maps remains invertible,
so the relative class is unchanged. The smoothed symbol is
$C^\infty$ in $(x,\xi)$, continuous in $\lambda$ in every symbol
seminorm, and exactly equal to the identity over $T^*U\times\Lambda$.

Homogeneous extension and quantization are performed with a single cover,
a single partition, and the same cutoffs for all parameters. In each
chart, Calderón--Vaillancourt estimates the difference of two
quantized operators by finitely many seminorms of the difference
of their symbols. Therefore,
\begin{equation}
 \lambda\longmapsto P_\lambda
 \quad\text{is continuous in}\quad
 \mathcal L\bigl(L^2(\mathbb R^n;\mathbb C^N)\bigr).
 \label{eq:continuidad-norma-cuantizacion-familiar}
\end{equation}
Choose a single
$\chi\in C_c^\infty(\mathbb R^n)$ equal to $1$ on the projection of the
common support, and replace the quantized family $A_\lambda$ by
\begin{equation}
 P_\lambda:=I+\chi(A_\lambda-I)\chi.
 \label{eq:localizacion-uniforme-familias}
\end{equation}
The Schwartz kernel of $P_\lambda-I$ is contained in
$\operatorname{supp}\chi\times\operatorname{supp}\chi$, uniformly
in $\lambda$. This proves realization.

We now prove descent. The same model by triples applies to the
product of the preceding compact Hausdorff pair with $[0,1]$. Thus a
relative homotopy of classes is represented, after stabilization,
over
\[
 B^*S^n\times\Lambda\times[0,1]
\]
and is elementary over $B^*U\times\Lambda\times[0,1]$. Repeat the
preceding smoothing on the compact space $\Lambda\times[0,1]$ and quantize
with the same cutoff $\chi$. This gives a norm-continuous family of
Fredholm operators. There exists a common finite-dimensional subspace $V$
making all stabilizations surjective, since
$\Lambda\times[0,1]$ is compact. Their kernels form a bundle
$\mathcal K\to\Lambda\times[0,1]$, and the index of the family is
\[
 [\mathcal K]-[\underline V]
 \in K^0(\Lambda\times[0,1]).
\]
The inclusions
$i_0,i_1:\Lambda\longrightarrow\Lambda\times[0,1]$ are homotopic;
consequently,
\[
 i_0^*\bigl([\mathcal K]-[\underline V]\bigr)
 =i_1^*\bigl([\mathcal K]-[\underline V]\bigr).
\]
Thus the index bundle is invariant under relative homotopies.

An isomorphism of triples depending on the covariable reduces, through
$(x,\xi,\lambda,t)\mapsto(x,t\xi,\lambda)$, to bundle isomorphisms over
the base; the procedure is uniform because $\Lambda$ is compact. An
elementary triple contracts in the same way to a family of
isomorphisms and has zero index bundle. Direct sum of triples corresponds
to direct sum of stabilized families. Thus the index respects equivalences, homotopies, elementary triples, and sums, which are all
the relations in the relative model. This proves that
\eqref{eq:indice-familiar-sobre-k} is well defined.
\end{proof}

The Bott isomorphism with parameters is
\[
 \alpha_{n,\Lambda}:
 K_c^0(\mathbb R^{2n}\times\Lambda)
 \longrightarrow K^0(\Lambda).
\]

\begin{theorem}[Euclidean index for families]
\label{teo:indice-euclidiano-familias}
Under the preceding hypotheses,
\begin{equation}
 \operatorname{Ind}_{\mathrm a}(P)
 =(-1)^n\alpha_{n,\Lambda}
 \bigl([\boldsymbol{\sigma}(P_\Lambda)]\bigr)
 \in K^0(\Lambda).
 \label{eq:formula-indice-euclidiano-familias}
\end{equation}
\end{theorem}

\begin{proof}
Proposition~\ref{prop:realizacion-descenso-familias-euclidianas}
constructs the homomorphism $I_{n,\Lambda}$ and proves its independence of
the relative representative and quantization.
This homomorphism is $K^0(\Lambda)$--linear. Indeed, if
$\mathbf{G}\longrightarrow\Lambda$ is a vector bundle, the stabilization defining
the index of $P\otimes \mathbf{G}$ is the preceding stabilization tensored with
$\mathbf{G}$; consequently,
\[
 \operatorname{Ind}_{\mathrm a}(P\otimes \mathbf{G})
 =\operatorname{Ind}_{\mathrm a}(P)\,[\mathbf{G}].
\]

By Theorem~\ref{teo:periodicidad-bott-k}, every class in the domain has
a unique expression
\[
 b^{\boxtimes n}\boxtimes a,
 \qquad a\in K^0(\Lambda).
\]
Iterating the angular regularization of
\eqref{eq:producto-regularizado-angular} and relative realization,
we obtain an operator
\[
 R_n\in\operatorname{Ell}_c^0(\mathbb R^n)
\]
such that
\[
 [\boldsymbol{\sigma}(R_n)]=b^{\boxtimes n},
 \qquad
 \operatorname{ind}R_n=(-1)^n.
\]

Write $a=[G]-[H]$ and represent $G$ and $H$, after adding trivial
bundles if necessary, by continuous projections
\[
 p_G:\Lambda\longrightarrow M_N(\mathbb C),
 \qquad
 p_H:\Lambda\longrightarrow M_N(\mathbb C).
\]
The family
\begin{equation}
 \mathcal R_G(\lambda)
 :=R_n\otimes p_G(\lambda)
   +I\otimes\bigl(I-p_G(\lambda)\bigr)
 \label{eq:familia-bott-proyeccion-G}
\end{equation}
is the identity on the complement of $G_\lambda$. Its symbol class and
index are
\[
 [\boldsymbol{\sigma}(\mathcal R_G)]
 =b^{\boxtimes n}\boxtimes[G],
 \qquad
 \operatorname{Ind}_{\mathrm a}(\mathcal R_G)=(-1)^n[G].
\]
For the negative term, use
\begin{equation}
 \mathcal R_H^-(\lambda)
 :=R_n^*\otimes p_H(\lambda)
   +I\otimes\bigl(I-p_H(\lambda)\bigr).
 \label{eq:familia-bott-proyeccion-H}
\end{equation}
Proposition~\ref{prop:adjunto-simbolo-k} and the equality
$\operatorname{ind}R_n^*=-\operatorname{ind}R_n$ give
\[
 [\boldsymbol{\sigma}(\mathcal R_H^-)]=-b^{\boxtimes n}\boxtimes[H],
 \qquad
 \operatorname{Ind}_{\mathrm a}(\mathcal R_H^-)=-(-1)^n[H].
\]
The direct sum of \eqref{eq:familia-bott-proyeccion-G} and
\eqref{eq:familia-bott-proyeccion-H} realizes
$b^{\boxtimes n}\boxtimes a$ and has index $(-1)^n a$. Therefore,
\[
 I_{n,\Lambda}(b^{\boxtimes n}\boxtimes a)=(-1)^n a.
\]
Since
\[
 \alpha_{n,\Lambda}(b^{\boxtimes n}\boxtimes a)=a,
\]
we obtain \eqref{eq:formula-indice-euclidiano-familias}.
\end{proof}

\section{Compatibility with embeddings of closed manifolds}
\label{sec:compatibilidades-indice-encaje}

The preceding constructions supply the three maps in the
embedding argument. If $j:X\hookrightarrow\mathbb R^N$ is an embedding of
a closed manifold and $a\in K_c^0(T^*X)$, we have
\[
 K_c^0(T^*X)
 \xrightarrow{\ \operatorname{Th}_{N(Tj)}\ }
 K_c^0(N(Tj))
 \xrightarrow{\ e\ }
 K_c^0(T^*\mathbb R^N)
 \xrightarrow{\ (-1)^N\alpha_N\ }
 \mathbb Z.
\]
The composition of the first two arrows is $j_!$, and the last is the
Euclidean index. Thus the topological expression associated with $a$ is
\begin{equation}
 \operatorname{ind}_{\mathrm t}(a)
 :=(-1)^N\alpha_N(j_!a).
 \label{eq:indice-topologico-preparacion-cerrada}
\end{equation}

Naturality of Thom, functoriality
\eqref{eq:functorialidad-pushforward-k}, and
\eqref{eq:inversa-pushforward-punto} show that
\eqref{eq:indice-topologico-preparacion-cerrada} is independent of the Euclidean
space in which the embedding is stabilized. Full independence of the
embedding and equality with the analytic index will be proved in
Chapter~\ref{cap:atiyah-singer-cerradas}. Analytically, the remaining step
is to realize the Thom class using the normal oscillator and to
prove that extension from a tubular neighborhood changes neither kernel nor
cokernel. Topologically, all the ingredients are already contained in
Theorems~\ref{teo:periodicidad-bott-k} and
\ref{teo:isomorfismo-thom-k} and in functoriality of the pushforward.

\chapter{The index theorem on closed manifolds}
\label{cap:atiyah-singer-cerradas}

Let $M$ be a closed smooth manifold of real dimension $n$, and let
$\mathbf{E}^+,\mathbf{E}^-\longrightarrow M$ be complex vector bundles. If
\[
P\colon \Gamma(M,\mathbf{E}^+)\longrightarrow \Gamma(M,\mathbf{E}^-)
\]
is a classical elliptic pseudodifferential operator,
Theorem~\ref{pseudo:teo-fredholm-eliptico} and
Corollary~\ref{pseudo:cor-independencia-s-indice} show that
\[
\operatorname{ind}_{\mathrm a}(P)
=\dim_{\mathbb C}\ker P-\dim_{\mathbb C}\ker P^*
\]
is well defined and independent of the Sobolev exponent used to
realize $P$ as a bounded operator.
Definition~\ref{def:clase-simbolo-k} associates a class to its principal symbol,
\[
[\boldsymbol{\sigma}(P)]\in K_c^0(T^*M).
\]
Theorem~\ref{teo:indice-analitico-sobre-k}, using the stability of
Theorem~\ref{pseudo:teo-estabilidad-indice}, proves that the index depends only
on this class. It remains to identify the resulting homomorphism
$K_c^0(T^*M)\longrightarrow\mathbb Z$.

The identification will be made using an embedding of $M$ into a
Euclidean space. The proof keeps two arguments separate. The first is
topological: the Thom isomorphism and extension by zero produce an
integer from a compactly supported class in $T^*M$. The second is
analytic: a Bott operator in the normal directions allows an
elliptic operator on $M$ to be transferred to a tubular neighborhood without changing its index.
The Euclidean formula of Chapter~\ref{cap:bott-thom-indice-euclidiano}
then completes the comparison. This is the embedding argument of
Atiyah and Singer; our reference exposition is
\cite[Chapter~12, pp.~284--301]{BleeckerBoossIndex}.

\begin{semblanzaHistorica}{Atiyah and Singer: two ways to count}
The theorem announced by Michael Atiyah and Isadore Singer in the 1960s states that two constructions of very different natures produce the same integer. The analytic index counts solutions and obstructions of an elliptic equation; the topological index combines the symbol class with operations in $K$--theory. The equality is more than a computational formula: it explains why global data of an equation remain invariant under deformation and why geometric operators recover invariants such as the Euler characteristic or characteristic genera~\cite{AtiyahSinger1968I}.
\end{semblanzaHistorica}

\section{The topological index}

Fix a Riemannian metric on $M$. The musical isomorphism
\[
\sharp_g\colon T^*M\longrightarrow TM,
\qquad
\xi\longmapsto \xi^{\sharp_g},
\]
induces an isomorphism
$K_c^0(T^*M)\cong K_c^0(TM)$. Although it will be used to write the construction,
the final result is independent of the metric: two metrics are joined by the
convex segment $g_t=(1-t)g_0+tg_1$, and the corresponding musical
isomorphisms form a proper homotopy of bundles.

By the Whitney embedding theorem in the form
\cite[Theorem~6.7(c), pp.~162--163]{BleeckerBoossIndex}, there exists a smooth
embedding
\[
j\colon M\hookrightarrow\mathbb R^N
\]
for some $N$. Let $\nu\longrightarrow M$ be its normal bundle, of rank
$r=N-n$. Theorem~\ref{teo:variedades-riemannianas-existencia-de-vecindades-tubulares}
provides open subsets $V\subseteq\nu$ and $U\subseteq\mathbb R^N$ containing,
respectively, the zero section and $j(M)$, and a diffeomorphism
\[
\Theta\colon V\longrightarrow U,
\qquad
\Theta(x,0)=j(x).
\]

The differential $Tj\colon TM\hookrightarrow T\mathbb R^N$ is also an
embedding. The differential of the projection of $\nu$ defines a real vector
bundle
\[
T\pi_\nu\colon T\nu\longrightarrow TM
\]
of rank $2r$: the fiber over $X\in T_xM$ consists of the vectors
tangent to $\nu$ that project to $X$. A normal connection separates the
horizontal and vertical directions and gives an identification
\[
T\nu\cong \pi_{TM}^*(\nu\oplus\nu)
\]
as bundles over $TM$. This is the normal bundle of $Tj(TM)$ in
$T\mathbb R^N$ under the differential of a tubular identification; in the
notation of Definition~\ref{def:pushforward-k-encaje},
$N(Tj)\cong(T\nu,T\pi_\nu)$.
The identification is not canonical, but two connections give
homotopic identifications. The complex structure
\begin{equation}
J(u,v)=(-v,u),
\qquad (u,v)\in\nu_x\oplus\nu_x,
\label{eq:estructura-compleja-normal-tm}
\end{equation},
however, determines a canonical complex orientation up to homotopy. By
Theorem~\ref{teo:isomorfismo-thom-k}, we obtain the Thom isomorphism
\[
\operatorname{Th}_{T\nu}\colon
K_c^0(TM)\longrightarrow K_c^0(T\nu).
\]
The tubular neighborhood identifies a neighborhood of $TM$ in
$T\mathbb R^N$ with a neighborhood of the zero section of $T\nu$. Extension
by zero then defines
\[
\operatorname{ext}\colon K_c^0(T\nu)\longrightarrow
K_c^0(T\mathbb R^N).
\]
For the structure \eqref{eq:estructura-compleja-normal-tm}, the complex
coordinate is $u+iv$ and the Thom class of a line is
$\lambda_{\mathbb C}=-b$ with the normalization of
Definition~\ref{def:clase-bott}. Thus, in rank $s$, we have
simultaneously
\begin{equation}
\lambda_{\mathbb C^s}=(-1)^s b^{\boxtimes s},
\qquad
\operatorname{ch}_c(\lambda_{\mathbf{W}})
=(-1)^s\Phi_{\mathbf{W}}\!\left(\operatorname{Td}(\overline{\mathbf{W}})^{-1}\right).
\label{eq:bloque-convenciones-thom}
\end{equation}
The complex orientation and Thom class are fixed by
\eqref{eq:bloque-convenciones-thom}.

In accordance with Definition~\ref{def:pushforward-k-encaje}, and using the
Euclidean metric to identify the last tangent bundle with the cotangent bundle,
the \emph{pushforward} of the embedding is
\begin{equation}
j_!:=\operatorname{ext}\circ\operatorname{Th}_{T\nu}\circ(\sharp_g)_*
\colon K_c^0(T^*M)\longrightarrow K_c^0(T^*\mathbb R^N).
\label{eq:pushforward-encaje-indice}
\end{equation}

Let $i_N\colon\{0\}\hookrightarrow\mathbb R^N$ be the inclusion. With the
conventions of Chapter~\ref{cap:bott-thom-indice-euclidiano}, if
$b\in K_c^0(\mathbb R^2)$ is the Bott class and
$\alpha_N\colon K_c^0(\mathbb R^{2N})\longrightarrow\mathbb Z$ is the
periodicity isomorphism, then
\begin{equation}
(i_N)_!(1)=(-1)^N b^{\boxtimes N},
\qquad
(i_N)_!^{-1}=(-1)^N\alpha_N.
\label{eq:normalizacion-pushforward-punto}
\end{equation}
These are precisely identities
\eqref{eq:pushforward-punto-bott} and
\eqref{eq:inversa-pushforward-punto}.
The sign comes from the complex structure
\eqref{eq:estructura-compleja-normal-tm}: in rank one, the Thom bundle homomorphism
is $z=v+i\eta$, whereas $b$ was normalized using the gluing function
$z^{-1}$.

\begin{definition}[Topological index]
\label{def:indice-topologico-cerradas}
\index{topological index}
Let $M$ be a closed smooth manifold and fix a smooth
embedding $j\colon M\hookrightarrow\mathbb R^N$. For
$a\in K_c^0(T^*M)$ define
\begin{equation}
\operatorname{ind}_{\mathrm t}(a)
:=(i_N)_!^{-1}(j_!a)
\in\mathbb Z.
\label{eq:def-indice-topologico-cerradas}
\end{equation}
\end{definition}

The definition appears to depend on the metric, embedding, and tubular
neighborhood. This dependence is only apparent.

\begin{proposition}[Independence of choices]
\label{prop:independencia-encaje-indice-topologico}
Let $M$ be a closed smooth manifold.
The integer in \eqref{eq:def-indice-topologico-cerradas} is independent of the
metric, normal connection, tubular neighborhood, and embedding
$j\colon M\hookrightarrow\mathbb R^N$. Moreover,
\[
\operatorname{ind}_{\mathrm t}\colon K_c^0(T^*M)\longrightarrow\mathbb Z
\]
is a homomorphism.
\end{proposition}

Independence under stabilization is stated in
\cite[Definition~12.11 and Exercise~12.10, p.~290]{BleeckerBoossIndex}; a
detailed proof appears in
\cite[Chapter~III, \S13, Definition~13.1]{LawsonMichelsohn1989}.

\begin{proof}
Independence of the metric was already obtained using the homotopy $g_t$.
The space of connections on a fixed bundle is affine; thus two normal
connections are joined by a segment and give homotopic Thom classes. Two
tubular neighborhoods admit, after restricting their domains, an
isotopy fixing the zero section. Homotopy invariance and excision
in compactly supported $K$--theory show that the extensions by zero
agree. These topological results are stated in
Theorem~\ref{teo:mapa-directo-k-apendice} of
Appendix~\ref{ap:topologia}.

It remains to compare two embeddings
$j_0\colon M\hookrightarrow\mathbb R^{N_0}$ and
$j_1\colon M\hookrightarrow\mathbb R^{N_1}$. Their stabilizations in
$\mathbb R^{N_0+N_1}$ are joined by the family
\[
j_\theta(x)=\bigl((\cos\theta)j_0(x),(\sin\theta)j_1(x)\bigr),
\qquad 0\leq\theta\leq\frac{\pi}{2}.
\]
Each $j_\theta$ is an embedding. At the endpoints this is immediate; if both
components are nonzero, injectivity of either component proves
injectivity, and the same observation applied to the differentials proves that
it is an immersion. Since $M$ is compact, an injective immersion into a
Hausdorff space is an embedding. Isotopy invariance gives equality
of the stabilized pushforwards. By functoriality of the Thom
isomorphism,
\[
(\iota\circ j)_!=\iota_!\circ j_!
\]
for a linear inclusion $\iota$ between Euclidean spaces. The same
identity holds for inclusions of points; hence the Bott factors
introduced by stabilization cancel when $(i_N)_!^{-1}$ is applied.
This proves independence of the embedding. Additivity follows from
additivity of the Thom isomorphism, extension by zero, and
$(i_N)_!^{-1}$.
\end{proof}

\section{Excision and products for the analytic index}

Passing from $M$ to a tubular neighborhood requires two analytic properties.
The first allows a symbol supported in an open subset to be transferred without changing
the index; the second calculates the index of an external product.

\begin{proposition}[Analytic excision]
\label{prop:excision-indice-analitico}
Let $X_0$ and $X_1$ be manifolds without boundary and let $U$ be an open set identified
with open subsets of both. For $i=0,1$, let $\mathbf{E}_i,\mathbf{F}_i\to X_i$ be complex vector bundles equipped with Hermitian bundle metrics
and identified over $U$, and let
$R_i\colon \mathbf{E}_i\to \mathbf{F}_i$ be a bundle homomorphism inducing a bounded
multiplication operator and corresponding, over $U$, to a common homomorphism
$R$. Suppose that there exist $L\Subset U$ such that $R_i$ is invertible outside
$L$ and a bounded operator
\[
B\colon L^2(U,\mathbf{E}|_U)\longrightarrow L^2(U,\mathbf{F}|_U)
\]
whose Schwartz kernel is supported in $L\times L$ and agrees under the
identifications. If
\[
P_i:=R_i+\operatorname{ext}_iB\colon
L^2(X_i,\mathbf{E}_i)\longrightarrow L^2(X_i,\mathbf{F}_i)
\]
is Fredholm for $i=0,1$, then
\begin{equation}
\operatorname{ind}P_0=\operatorname{ind}P_1.
\label{eq:excision-indice-analitico}
\end{equation}
In particular, the analytic index of a class in
$K_c^0(T^*U)$ extended by zero depends only on a neighborhood of the
projection of its support.
\end{proposition}

\begin{proof}
Choose on $U$ a Riemannian metric, its associated measure, and Hermitian
metrics on the bundles, and extend the metrics to $X_0$ and $X_1$. Thus the identifications over $U$ preserve
the inner products of $L^2(U,\mathbf E_0\restriction_U)$ and
$L^2(U,\mathbf F_0\restriction_U)$, and the adjoints. If $P_i u=0$, then
$R_i u=0$ outside $L$; invertibility of $R_i$ there implies
$\operatorname{supp}u\subseteq L$. The same argument for $P_i^*$ shows
$\operatorname{supp}v\subseteq L$ for every $v\in\ker P_i^*$.

Over $U$, the operators $P_0$ and $P_1$ agree with $R+B$, and the
same holds for their adjoints. Restriction therefore identifies
$\ker P_0$ with $\ker P_1$ and $\ker P_0^*$ with $\ker P_1^*$. Subtracting their
dimensions proves \eqref{eq:excision-indice-analitico}.

For a class $a\in K_c^0(T^*U)$, the extension-by-zero construction
allows us, after stabilization, to choose a localized representative of this form:
if $A$ is a quantization and $R$ its invertible part outside a compact set,
take $Q:=R+\chi(A-R)\chi$ with $\chi\in C_c^\infty(U)$ equal to one near
the support. The operator $B:=Q-R$ has kernel supported in a compact subset of
$U\times U$, and the preceding argument applies to its two extensions.
\end{proof}

To state the second property, it is convenient to use graded operators. If
$P\colon H_P^+\longrightarrow H_P^-$ and
$Q\colon H_Q^+\longrightarrow H_Q^-$ are closed densely defined
operators, write
\[
D_P=
\begin{pmatrix}
0&P^*\\
P&0
\end{pmatrix},
\qquad
D_Q=
\begin{pmatrix}
0&Q^*\\
Q&0
\end{pmatrix},
\qquad
\Gamma_P=
\begin{pmatrix}
I&0\\
0&-I
\end{pmatrix}.
\]
The odd operator
\begin{equation}
D_P\widehat\otimes I+\Gamma_P\widehat\otimes D_Q
\label{eq:producto-operadores-graduados}
\end{equation}
acts on the graded tensor product. Its component from even to odd
degree is
\begin{equation}
P\#Q=
\begin{pmatrix}
P\otimes I&-I\otimes Q^*\\
I\otimes Q&P^*\otimes I
\end{pmatrix}.
\label{eq:operador-producto-indice}
\end{equation}

\begin{theorem}[Analytic multiplicativity]
\label{teo:multiplicatividad-indice-analitico}
Let
$P\in\Psi^1_{\mathrm{cl}}(M;\mathbf{E}^+,\mathbf{E}^-)$ and
$Q\in\Psi^1_{\mathrm{cl}}(N;\mathbf{F}^+,\mathbf{F}^-)$ be elliptic operators on closed smooth
manifolds $M$ and $N$, respectively. The closure of the operator $P\#Q$ in
\eqref{eq:operador-producto-indice} is Fredholm and
\begin{equation}
\operatorname{ind}_{\mathrm a}(P\#Q)
=\operatorname{ind}_{\mathrm a}(P)
 \operatorname{ind}_{\mathrm a}(Q).
\label{eq:multiplicatividad-indice-analitico}
\end{equation}
\end{theorem}

\begin{proof}
Set
\[
 \mathcal H_P:=L^2(M,\mathbf{E}^+)\oplus L^2(M,\mathbf{E}^-),
 \qquad
 \mathcal H_Q:=L^2(N,\mathbf{F}^+)\oplus L^2(N,\mathbf{F}^-),
 \qquad
 \mathcal H:=\mathcal H_P\widehat\otimes\mathcal H_Q,
\]
with their natural orthogonal sums and gradings.
The anticommutation relation
$D_P\Gamma_P+\Gamma_PD_P=0$ implies
\begin{equation}
\left(D_P\widehat\otimes I+
\Gamma_P\widehat\otimes D_Q\right)^2
=D_P^2\widehat\otimes I+I\widehat\otimes D_Q^2.
\label{eq:cuadrado-producto-graduado}
\end{equation}
For every vector $u$ in the algebraic tensor product of smooth sections,
\[
\begin{aligned}
&\left\langle
\left(D_P^2\widehat\otimes I+I\widehat\otimes D_Q^2\right)u,u
\right\rangle_{\mathcal H}\\
&\qquad=\left\|(D_P\widehat\otimes I)u\right\|_{\mathcal H}^2
+\left\|(I\widehat\otimes D_Q)u\right\|_{\mathcal H}^2.
\end{aligned}
\]
Let us specify the closure used in these identities. Elliptic regularity
and Rellich imply that $D_P$ and $D_Q$, with domains
$H^1(M,\mathbf E^+\oplus\mathbf E^-)$ and
$H^1(N,\mathbf F^+\oplus\mathbf F^-)$, respectively, are self-adjoint
with compact resolvent. Let
$(\lambda_j)_{j\in I}$ and $(\mu_k)_{k\in J}$ be the distinct eigenvalues
of $D_P^2$ and $D_Q^2$, respectively, and let $V_j$ and $W_k$ be their
eigenspaces. The index sets $I,J$ are finite or countable.
Every eigenspace is finite-dimensional and invariant under the corresponding
operator and its grading. On $V_j\widehat\otimes W_k$, the
operator in \eqref{eq:producto-operadores-graduados} is a self-adjoint
matrix whose square is $(\lambda_j+\mu_k)I$.

The orthogonal sum of these matrices is self-adjoint on the domain
\[
\left\{\mathbf{u}=\sum_{(j,k)\in I\times J}\mathbf{u}_{jk}
\ \middle|\
\sum_{(j,k)\in I\times J}
(\lambda_j+\mu_k)\|\mathbf{u}_{jk}\|^2<\infty\right\},
\]
where the first sum converges in $\mathcal H$. Vectors with finitely
many components form a core in the graph norm, by
truncation of the second sum. This proves that the orthogonal sum is the
closure of the product initially defined on the algebraic tensor
product of smooth sections. For every $R>0$, there are only finitely
many pairs $(j,k)$ with $\lambda_j+\mu_k\leq R$, since both spectral
sequences have this property. Zero therefore has finite
multiplicity and is separated from the rest of the spectrum. The product is Fredholm.
The only blocks contributing to its kernel are those satisfying
$\lambda_j=\mu_k=0$. Hence
\[
\ker\!\left(D_P\widehat\otimes I+
\Gamma_P\widehat\otimes D_Q\right)
=\ker D_P\widehat\otimes\ker D_Q\subset\mathcal H.
\]
Separating even and odd components gives
\begin{align*}
\ker(P\#Q)
&\cong
(\ker P\otimes\ker Q)
\oplus(\ker P^*\otimes\ker Q^*),\\
\ker(P\#Q)^*
&\cong
(\ker P^*\otimes\ker Q)
\oplus(\ker P\otimes\ker Q^*).
\end{align*}
Subtracting dimensions yields
\[
\bigl(\dim\ker P-\dim\ker P^*\bigr)
\bigl(\dim\ker Q-\dim\ker Q^*\bigr),
\]
which is \eqref{eq:multiplicatividad-indice-analitico}.
\end{proof}

\section{The normal Bott operator}

The preceding multiplicativity must be applied uniformly in the
fibers of the normal bundle. We begin with the model on $\mathbb R^r$. On
$\Lambda^*\mathbb C^r$ let $\varepsilon_j$ be exterior multiplication by
the $j$th basis vector and $\iota_j=\varepsilon_j^*$ contraction.
On the Schwartz class define
\begin{equation}
\mathcal B_r
=\sum_{j=1}^r
\left(
\varepsilon_j\left(\frac{\partial}{\partial v_j}+v_j\right)
+\iota_j\left(-\frac{\partial}{\partial v_j}+v_j\right)
\right).
\label{eq:operador-bott-oscilador}
\end{equation}
This is an odd symmetric operator. Denote by
$\mathcal B_r^+$ its component from even to odd forms.

\begin{lemma}[Kernel and index of the Bott operator]
\label{lem:nucleo-operador-bott-normal}
The closure of $\mathcal B_r$ is self-adjoint, has compact resolvent, and
\[
\ker\mathcal B_r
=\mathbb C e^{-\frac{|v|^2}{2}}\otimes\Lambda^0\mathbb C^r.
\]
In particular,
\[
\operatorname{ind}_{\mathrm a}(\mathcal B_r^+)=1.
\]
Its total symbol in the variables $(v,\eta)\in\mathbb R^r\oplus\mathbb R^r$
is
\begin{equation}
\beta_r(v,\eta)
=\varepsilon(v+i\eta)+\varepsilon(v+i\eta)^*,
\label{eq:simbolo-total-bott-normal}
\end{equation}
where $\displaystyle\varepsilon(z)=\displaystyle\sum_{j=1}^r z_j\varepsilon_j$ and
$\displaystyle\varepsilon(z)^*=\displaystyle\sum_{j=1}^r\overline{z_j}\iota_j$.
Its component from even to odd degree represents the Thom class of
$\mathbb R^r\oplus\mathbb R^r$ with the complex structure
$J(v,\eta)=(-\eta,v)$.
\end{lemma}

\begin{proof}
The relations
\[
\varepsilon_i\varepsilon_j+\varepsilon_j\varepsilon_i=0,
\qquad
\iota_i\iota_j+\iota_j\iota_i=0,
\qquad
\varepsilon_i\iota_j+\iota_j\varepsilon_i=\delta_{ij}I
\]
give, on $\mathcal S(\mathbb R^r,\Lambda^*\mathbb C^r)$,
\begin{equation}
\mathcal B_r^2
=-\Delta_v+|v|^2-r+2\mathcal N,
\qquad
\mathcal N=\sum_{j=1}^r\varepsilon_j\iota_j.
\label{eq:cuadrado-operador-bott}
\end{equation}
The operator $\mathcal N$ acts by multiplication by $q$ on
$\Lambda^q\mathbb C^r$. The scalar oscillator
$-\Delta_v+|v|^2$ has eigenvalues
$r+2|\alpha|$, $\alpha\in\mathbb N_0^r$, and its ground state is
$e^{-\frac{|v|^2}{2}}$. Thus the eigenvalues of
\eqref{eq:cuadrado-operador-bott} are $2|\alpha|+2q$. Zero occurs exactly
once, for $\alpha=0$ and $q=0$. The Hermite function basis proves
both essential self-adjointness and compactness of the resolvent.

The first-order part together with the linear potential gives
\eqref{eq:simbolo-total-bott-normal}. Since
\[
\beta_r(v,\eta)^2=(|v|^2+|\eta|^2)I,
\]
its restriction to the sphere is an isomorphism between the even and
odd parts. To compare it with the sign in
Definition~\ref{def:clase-thom-compleja}, set $z=v+i\eta$ and consider
$c_t(z)=\varepsilon(z)+e^{i\pi t}\varepsilon(z)^*$,
$0\leq t\leq1$. The creation and annihilation relations give
$c_t(z)^2=e^{i\pi t}|z|^2I$. Thus the even-to-odd components are
invertible when $z\neq0$ and join the preceding symbol to
$\varepsilon(z)-\varepsilon(z)^*$, which defines the Thom class. The
normalization agrees with \eqref{eq:normalizacion-pushforward-punto}: for
$r=1$, the bundle homomorphism is $v+i\eta$ and represents $-b$.
\end{proof}

Formula \eqref{eq:operador-bott-oscilador} is invariant under the natural
action of $O(r)$. Indeed, differentiation, multiplication by $v$, and the
exterior and interior operations transform simultaneously. This
equivariance allows us to pass from the vector-space model to the normal bundle without
choosing global frames.

Existence of a global representative on the fiber compactification
uses the equivariant Bott normalization for the standard representation of $O(r)$, stated in
Theorem~\ref{teo:bott-equivariante-apendice}. The result provides
the class and its index in $R(O(r))$; the next lemma develops the
analytic normalization we need and adapts it to the normal bundle. Before
stating it, let us specify the compactification used in this step.

Let $G=O(r)$ and let $S^r=\mathbb R^r\cup\{\infty\}$, with the action of $G$
fixing $\infty$. Denote the orthonormal frame bundle of $\nu$ by
$\boldsymbol{\mathcal P}_\nu\longrightarrow M$. The fiberwise spherical
compactification is the closed manifold
\begin{equation}
\widehat\nu:=\boldsymbol{\mathcal P}_\nu\times_G S^r,
\qquad
\widehat\pi:\widehat\nu\longrightarrow M.
\label{eq:compactificacion-esferica-normal}
\end{equation}
The section $\mathbf{s}_\infty=\boldsymbol{\mathcal P}_\nu\times_G\{\infty\}$ is the complement
of the open copy of $\nu$ in $\widehat\nu$.

\begin{lemma}[Equivariant normal representative]
\label{lem:representante-normal-equivariante}
Let $S^r$ be the closed smooth sphere and let $G=O(r)$ with the
action described above.
There exist complex bundles $\mathbf{W}^+,\mathbf{W}^-\longrightarrow S^r$ with an action of $G$ and
a classical elliptic $G$--equivariant pseudodifferential operator of order one
\[
B\colon\Gamma(S^r,\mathbf{W}^+)\longrightarrow\Gamma(S^r,\mathbf{W}^-)
\]
such that its symbol class is the extension by zero from
$T^*\mathbb R^r$ of the class of \eqref{eq:simbolo-total-bott-normal}.
The equality holds in equivariant $K$--theory. Moreover, after modifying $B$ by
an $G$--equivariant smoothing operator, we may assume that
\begin{equation}
\ker B=\mathbb C,
\qquad
\ker B^*=0,
\label{eq:nucleos-bott-compactificado}
\end{equation}
where $G$ acts trivially on $\mathbb C$.
\end{lemma}

\begin{proof}
The bounded transform
\[
F_r^+=\mathcal B_r^+
\bigl(I+\mathcal B_r^-\mathcal B_r^+\bigr)^{-\frac{1}{2}}
\]
is Fredholm and $G$--equivariant. The calculation in
Lemma~\ref{lem:nucleo-operador-bott-normal} shows, as representations of
$G$ and not merely after taking dimensions, that
\begin{equation}
[\ker F_r^+]-[\ker(F_r^+)^*]=[\mathbb C]
\quad\text{in }R(G).
\label{eq:indice-equivariante-bott}
\end{equation}
To pass to the compactification, apply
Theorem~\ref{teo:bott-equivariante-apendice} of
Appendix~\ref{ap:topologia}. If $i:\{0\}\hookrightarrow\mathbb R^r$, this
result provides an elliptic $G$--equivariant representative of
$\beta_r=i_!(1)$ on $S^r$, elementary near $\infty$, and establishes
\begin{equation}
\operatorname{Ind}_G(\beta_r)=[\mathbb C]
\quad\text{in }R(G),
\label{eq:lema-bott-equivariante-indice}
\end{equation}
with the trivial representation on the right-hand side.

The normalization in \eqref{eq:lema-bott-equivariante-indice} agrees with the one
adopted here. Indeed, the local representative is
\eqref{eq:simbolo-total-bott-normal}, and
\eqref{eq:indice-equivariante-bott} calculates its index as a representation,
not just its dimension: the Gaussian is fixed by $O(r)$ and the cokernel is
zero. This rules out the determinant representation and fixes the right-hand
side of \eqref{eq:lema-bott-equivariante-indice}.

Let $B_0$ be an order-zero quantization of the symbol provided by
that theorem. We can make it $G$--equivariant by averaging its conjugates
with respect to normalized Haar measure on $G$: the principal symbol
remains unchanged because it is equivariant. Then
\[
[\ker B_0]-[\ker B_0^*]=[\mathbb C]
\quad\text{in }R(G).
\]
Since $G$ is compact, its
finite-dimensional representations are semisimple. The preceding equality
provides a $G$--equivariant
decomposition
\[
\ker B_0\cong\mathbb C\oplus\ker B_0^*.
\]
Let $K$ be the operator identifying the second summand with $\ker B_0^*$
and vanishing on its orthogonal complement. It has finite rank and is
$G$--equivariant. Elliptic regularity of both kernels shows that
its integral kernel is smooth. The operator $B_1:=B_0+K$ is surjective:
$B_0$ identifies $(\ker B_0)^\perp$ with its closed image, and $K$ covers
the orthogonal complement of that image. Its kernel is the first copy
of $\mathbb C$. Finally, compose on the right with the positive,
invertible, $G$--equivariant operator
\[
\Lambda:=
\bigl(I+(\nabla^{\mathbf{W}^+})^*\nabla^{\mathbf{W}^+}\bigr)^{1/2},
\]
of order one, using the round metric and a $G$--invariant Hermitian metric
and connection on $\mathbf{W}^+$. The map
$\Lambda:H^1(S^r,\mathbf{W}^+)\to L^2(S^r,\mathbf{W}^+)$ is an
isomorphism by the calculus of powers in
Chapter~\ref{cap:operadores-pseudodiferenciales}. Thus
$B:=B_1\Lambda$ has the same image as $B_1$ and
$\ker B=\Lambda^{-1}(\ker B_1)$. This identification is
$G$--equivariant and preserves the one-dimensional trivial kernel. The positive
principal symbol of $\Lambda$ does not change the symbol class after
normalizing the order. This gives the operator $B$ in the statement.
\end{proof}

Fix a metric $h_\nu$ on $\nu$. A metric connection on $\nu$ splits the cotangent bundle of
$\widehat\nu$ into horizontal and vertical parts,
\begin{equation}
T^*\widehat\nu=\mathbf{H}^*\oplus \mathbf{V}^*,
\qquad
\mathbf{H}^*\cong\widehat\pi^*T^*M,
\qquad
\mathbf{V}^*\cong\boldsymbol{\mathcal P}_\nu\times_G T^*S^r.
\label{eq:descomposicion-cotangente-esfera-normal}
\end{equation}
Denote by $|\cdot|_{\mathbf{H}}$ and $|\cdot|_{\mathbf{V}}$ the norms induced on
$\mathbf{H}^*$ and $\mathbf{V}^*$, and by $|\cdot|_{\widehat g}$ the norm of the orthogonal sum.
The operator $B$ induces a vertical operator $B_\nu$ on
$\widehat\nu$, between the associated bundles
\[
\mathbf{W}_\nu^\pm:=\boldsymbol{\mathcal P}_\nu\times_G \mathbf{W}^\pm
\longrightarrow\widehat\nu.
\]
Its kernel is the trivial line over $M$ and its cokernel is zero, by
\eqref{eq:nucleos-bott-compactificado}.

The product construction will use cutoffs depending on the vertical
operator. These cutoffs have two properties that we will need simultaneously:
they preserve its spectral subspaces and turn local products into
classical operators on the total manifold. We begin by specifying the
domain corresponding to the latter.

\begin{lemma}[Domain of a classical elliptic operator]
\label{lem:estimacion-conica-horizontal-vertical-parametro}
Let $X$ be a closed smooth manifold and let $\mathbf{E},\mathbf{F}\to X$
be complex vector bundles equipped with Hermitian bundle metrics. If
$A\in\Psi^1_{\mathrm{cl}}(X;\mathbf{E},\mathbf{F})$ is elliptic, its closure
in $L^2(X,\mathbf E)$ has domain $H^1(X,\mathbf{E})$, and the domain of its Hilbert space
adjoint is $H^1(X,\mathbf{F})$. There exist constants $C_0,C_1>0$ such
that, for every $\mathbf{u}\in H^1(X,\mathbf{E})$,
\begin{align}
\|A\mathbf{u}\|_{L^2(X,\mathbf{F})}
&\leq C_0\|\mathbf{u}\|_{H^1(X,\mathbf{E})},
\label{eq:cota-superior-conica-parametro}\\
\|\mathbf{u}\|_{H^1(X,\mathbf{E})}
&\leq C_1\bigl(\|A\mathbf{u}\|_{L^2(X,\mathbf{F})}
+\|\mathbf{u}\|_{L^2(X,\mathbf{E})}\bigr).
\label{eq:cota-eliptica-conica-parametro}
\end{align}
The corresponding assertions hold for the formal adjoint. In
particular, a formally self-adjoint classical elliptic operator
of order one is self-adjoint with domain $H^1(X,\mathbf E)$.
\end{lemma}

\begin{proof}
The first inequality is
Corollary~\ref{pseudo:cor-mapeo-sobolev-cerrada}. A parametrix
$Q\in\Psi^{-1}_{\mathrm{cl}}(X;\mathbf{F},\mathbf{E})$ satisfies
$QA=I-K$, where $K$ is smoothing. Applying the identity
$\mathbf{u}=QA\mathbf{u}+K\mathbf{u}$ and the same corollary yields
\eqref{eq:cota-eliptica-conica-parametro}. If $\mathbf u_j\in\Gamma(\mathbf E)$ converges to $\mathbf u$ in
$H^1(X,\mathbf E)$, the first bound gives
$A\mathbf u_j\to A\mathbf u$ in $L^2(X,\mathbf F)$; thus
$\mathbf u$ belongs to the domain of the closure.
Conversely, if $(\mathbf u_j,A\mathbf u_j)$ converges in
$L^2(X,\mathbf E)\times L^2(X,\mathbf F)$ to $(\mathbf u,\mathbf f)$,
the second bound applied to $\mathbf u_j-\mathbf u_k$ proves that
$(\mathbf u_j)$ is Cauchy in $H^1(X,\mathbf E)$.
Its limit in that space agrees with $\mathbf u$ by the continuous
inclusion into $L^2(X,\mathbf E)$; the first bound identifies
$A\mathbf u=\mathbf f$. This proves both inclusions for the domain.

If $\mathbf{v}$ belongs to the domain of the Hilbert space adjoint, the identity
defining that adjoint says that $A_h^*\mathbf{v}$ belongs to $L^2(X,\mathbf E)$ in the
distributional sense. The parametrix of the classical elliptic operator $A_h^*$
then gives $\mathbf{v}\in H^1(X,\mathbf{F})$. Conversely, let $\mathbf v\in H^1(X,\mathbf F)$ and take
$\mathbf v_j\in\Gamma(\mathbf F)$ with
$\mathbf v_j\to\mathbf v$ in $H^1(X,\mathbf F)$.
The mapping bound for $A_h^*$ gives
$A_h^*\mathbf v_j\to A_h^*\mathbf v$ in $L^2(X,\mathbf E)$.
For each $\mathbf u\in\Gamma(\mathbf E)$, Cauchy--Schwarz allows us to
pass to the limit in
\[
 \langle A\mathbf u,\mathbf v_j\rangle_{L^2(X,\mathbf F)}
 =\langle\mathbf u,A_h^*\mathbf v_j\rangle_{L^2(X,\mathbf E)}.
\]
The resulting identity and
$|\langle\mathbf u,A_h^*\mathbf v\rangle_{L^2(X,\mathbf E)}|
\leq\|\mathbf u\|_{L^2(X,\mathbf E)}\|A_h^*\mathbf v\|_{L^2(X,\mathbf E)}$
are precisely the definition of the domain of the Hilbert space adjoint. The two operators agree there. If $A_h^*=A$, their
domains and actions agree, which is the asserted self-adjointness.
\end{proof}

\begin{lemma}[Compactness from the first-order domain]
\label{lem:compacidad-orden-cero-h1-l2}
Let $X$ be a closed smooth manifold and let $\mathbf{E},\mathbf{F}\to X$ be complex vector bundles
equipped with Hermitian bundle metrics. If $T:L^2(X,\mathbf{E})\to L^2(X,\mathbf{F})$ is bounded,
its restriction is compact as a map
\[
T:H^1(X,\mathbf{E})\longrightarrow L^2(X,\mathbf{F}).
\]
This applies, in particular, to classical operators of order zero.
\end{lemma}

\begin{proof}
The restriction is the composition of the compact inclusion
$H^1(X,\mathbf{E})\hookrightarrow L^2(X,\mathbf{E})$, given by
Theorem~\ref{teo: rellich--kondrashov para haces vectoriales}, with the
bounded operator $T$. For classical operators of order zero, continuity
from $L^2(X,\mathbf E)$ to $L^2(X,\mathbf F)$ is Corollary~\ref{pseudo:cor-mapeo-sobolev-cerrada}.
\end{proof}

\begin{lemma}[Spectral cutoffs on a local product]
\label{lem:cortes-espectrales-producto-clasico}
Let $n\in\mathbb N$, let $Z$ be a closed smooth manifold, and let $\mathbf{W}\to Z$ be a complex vector bundle
equipped with a Hermitian bundle metric. Suppose that
$D\in\Psi^1_{\mathrm{cl}}(Z;\mathbf{W},\mathbf{W})$ is elliptic and
self-adjoint with domain $H^1(Z,\mathbf{W})$. Set $L=I+D^2$ and
denote the principal symbol of $D$ by $b_1(z,\eta)$. Fix
$h\in C_c^\infty(\mathbb R)$ such that
\[
0\leq h\leq1,\qquad
h(t)=1\quad(-1\leq t\leq1),\qquad
h(t)=0\quad(t\geq4).
\]
For $R\geq1$ define $\Phi_R$ using the partial Fourier transform
in $x\in\mathbb R^n$:
\begin{equation}
\widehat{\Phi_R\mathbf{u}}(\xi)
=h\!\left(\frac{R^2L}{1+|\xi|^2}\right)\widehat{\mathbf{u}}(\xi),
\qquad \xi\in\mathbb R^n.
\label{eq:corte-conico-producto}
\end{equation}
It is self-adjoint and $0\leq\Phi_R\leq I$ on
$L^2(\mathbb R^n\times Z,\mathbf{W})$. After multiplication on the left
and right by smooth compactly supported functions in $x$, the operators
\[
\Phi_R,\qquad
(P\otimes I)\Phi_R,\qquad
\Phi_R(P\otimes I),\qquad
(I\otimes D)(I-\Phi_R)
\]
are classical in the joint variables $(x,z)$, of orders $0,1,1,1$,
respectively, for every classical operator $P$ of order one in
$\mathbb R^n$ between trivial bundles of finite rank. Identities on the additional factors
are understood in the formulas.

The principal symbol of the cutoff is
\begin{equation}
f_R(z;\xi,\eta)
=h\!\left(\frac{R^2b_1(z,\eta)^2}{|\xi|^2}\right)
\quad(\xi\neq0),
\label{eq:simbolo-principal-corte-espectral}
\end{equation}
extended by zero when $\xi=0$ and $\eta\neq0$. Near the axis
$\eta=0$, with $\xi\neq0$, it is the identity. The other three principal
symbols are $p_1f_R$, $f_Rp_1$, and $b_1(I-f_R)$.

Every bounded operator in the vertical variable commuting with the resolvent
of $D$ commutes with $\Phi_R$. Every grading
anticommuting with $D$ also commutes with $\Phi_R$. If a compact group acts unitarily
on $\mathbf{W}$ and $D$ is equivariant, the cutoffs are equivariant.
\end{lemma}

\begin{proof}
By the spectral theorem, for every $\xi\in\mathbb R^n$ the operator
$h(R^2L/(1+|\xi|^2))$ is a positive self-adjoint contraction. Plancherel
in the variable $x$ proves the same properties for $\Phi_R$. The
commutation relations follow from functional calculus; if $\Gamma D=-D\Gamma$,
then $\Gamma$ commutes with $L$ and all its spectral functions.

Let us prove the assertion about joint symbols. We use the construction
of the resolvent parametrix in
Proposition~\ref{pseudo:prop-parametrix-resolvente}, retaining the
dependence of the bounds on distance to the real axis. Write
\[
\tau=R^{-1}(1+|\xi|^2)^{1/2},\qquad
h\!\left(\frac{R^2L}{1+|\xi|^2}\right)=h(\tau^{-2}L).
\]
In a chart of $Z$, the symbol of $L$ has an expansion
$\ell\sim\displaystyle\sum_{j=0}^{\infty}\ell_{2-j}$, whose first
term is $\ell_2=b_1^2$. It is a positive Hermitian matrix for
$\eta\neq0$. Compactness and ellipticity give constants
$c_Z,C_Z>0$ such that
\begin{equation}
c_Z|\eta|^2I\leq\ell_2(z,\eta)\leq C_Z|\eta|^2I.
\label{eq:cotas-principal-vertical-cuadrado}
\end{equation}

To integrate the resolvent against a smooth function, fix
$N\in\mathbb N$ and a cutoff $\vartheta\in C_c^\infty(\mathbb R)$ equal
to one near zero. With $w=a+it$ set
\[
\widetilde h_N(a+it)
=\vartheta(t)\sum_{j=0}^{N}\frac{(it)^j}{j!}h^{(j)}(a).
\]
The operator $\partial_{\overline w}=\frac12(\partial_a+i\partial_t)$
satisfies
$|\partial_{\overline w}\widetilde h_N(a+it)|\leq C_N|t|^N$ for
$|t|$ sufficiently small. The Cauchy--Green formula, first for
a real number and then integrated against the spectral resolution of
$L$, gives
\begin{equation}
h(\tau^{-2}L)
=\frac1\pi\int_{\mathbb R^2}
\partial_{\overline w}\widetilde h_N(w)\,
\tau^2(L-\tau^2w)^{-1}\,da\,dt.
\label{eq:cauchy-green-corte-vertical}
\end{equation}
The integrand is considered for $t\neq0$; the set $t=0$ has measure
zero. The spectral bound
\[
\|\tau^2(L-\tau^2w)^{-1}\|_{\mathcal L(L^2(Z,\mathbf{W}))}
\leq |t|^{-1}
\]
justifies the integral in operator norm when $N\geq1$.

The symbol recurrence of
Chapter~\ref{cap:operadores-pseudodiferenciales} now applies to
$L-\tau^2w$. In the region $\eta\neq0$ it starts with
\[
r_{-2}=(\ell_2-\tau^2w)^{-1}
\]
and, for $j\in\mathbb N$, continues by
\begin{equation}
r_{-2-j}
=-r_{-2}
\sum_{\substack{a,b\in\mathbb N_0,\ \gamma\in\mathbb N_0^{\dim Z}\\
 a+b+|\gamma|=j,\ b<j}}
\frac{i^{-|\gamma|}}{\gamma!}
(\partial_\eta^\gamma\ell_{2-a})
(\partial_z^\gamma r_{-2-b}).
\label{eq:recurrencia-corte-espectral-vertical}
\end{equation}
On an annulus where $|\eta|$ and $\tau$ are comparable, all
derivatives of $r_{-2-j}$ satisfy the parameter-dependent estimates of order
$-2-j$. Their constants are bounded by a power of $|t|^{-1}$:
this is seen first for $r_{-2}$ from
\eqref{eq:cotas-principal-vertical-cuadrado}, and when the inverse is differentiated,
each new derivative adds finitely many factors $r_{-2}$. The finite
sum in \eqref{eq:recurrencia-corte-espectral-vertical} retains this
property. Specifically, for every finite set of derivatives, there exists an
integer $m_N$ bounding the power of $|t|^{-1}$ needed both in the
finite parametrix and in its remainder.

In \eqref{eq:cauchy-green-corte-vertical} choose an extension of order
greater than $m_N+2$. We can then integrate each term and its remainder.
The left and right parametrix identities, followed by the bound for
the exact resolvent, control the difference between the two resolvents:
if $(L-\tau^2w)Q_N=I-S_N$, then
\[
(L-\tau^2w)^{-1}-Q_N=(L-\tau^2w)^{-1}S_N.
\]
Choose the expansion length so that $S_N$ has the required negative
order before integration. This produces, on the annulus
$|\eta|\asymp\tau$, a classical expansion whose first term is
$h(\ell_2/\tau^2)$.

Let us also examine what happens outside this annulus. Every term in the
recurrence is a rational function of $w$; its poles lie in the spectrum
of $\ell_2/\tau^2$. If $|\eta|/\tau$ is sufficiently small, this
spectrum lies in an interval where $h=1$. The integral of the first term
is then $I$. Those of the subsequent terms vanish: they are integrals
of products of resolvents whose expansion in $w$ starts with $w^{-2}$,
and integrating a constant function around all its poles gives
zero. This argument remains valid for noncommuting matrices,
since it applies entry by entry to the complete rational function. If
$|\eta|/\tau$ is sufficiently large, $h$ vanishes near all the
poles and all the integrals vanish. Introducing cutoffs in these three
regions leaves the terms differentiating the cutoffs in the transition
annuli.

To control the remainders near the axis $\eta=0$ as well, separate the
frequencies $|\eta|\leq2\tau^{1/2}$, assuming first that $\tau\geq1$.
Let $Q_\tau$ be a vertical localization whose symbol is one for
$|\eta|\leq\tau^{1/2}$ and zero for $|\eta|\geq2\tau^{1/2}$, with
the spatial cutoffs of the chart. For $N\in\mathbb N$ define
\[
g_N(t)=\frac{h(t)-1}{t^N}\quad(t>0),\qquad g_N(0)=0.
\]
This function is smooth and bounded on $[0,+\infty)$, since its numerator
vanishes near zero. Spectral calculus gives the exact identity
\[
\bigl(h(\tau^{-2}L)-I\bigr)Q_\tau
=\tau^{-2N}g_N(\tau^{-2}L)L^NQ_\tau.
\]
Vertical Sobolev norms can be expressed using powers of
$L$, by Proposition~\ref{pseudo:prop-dominios-potencias-complejas}.
Thus $g_N(\tau^{-2}L)$ is uniformly bounded on each of
these spaces. For $s\in\mathbb N_0$, symbol calculus and the
bound $\langle\eta\rangle\leq C\tau^{1/2}$ on the cutoff support
give
\[
\|L^NQ_\tau\|_{\mathcal L(H^{-s}(Z,\mathbf{W}),H^s(Z,\mathbf{W}))}
\leq C_{N,s}\tau^{N+s}.
\]
Indeed, conjugating by the powers of $L$ defining the two norms
produces a symbol of order $2N+2s$ with this support; dividing it by
$\tau^{N+s}$ leaves bounded the order-zero seminorms needed for
the estimate on $L^2(Z,\mathbf W)$. We thus obtain
\[
\|\bigl(h(\tau^{-2}L)-I\bigr)Q_\tau\|_{
\mathcal L(H^{-s}(Z,\mathbf{W}),H^s(Z,\mathbf{W}))}
\leq C_{N,s}\tau^{-N+s}.
\]
The integer $N$ is arbitrary. Differentiating in $\tau$ produces derivatives of the
cutoff and functions obtained from $g_N$ by applying $t\partial_t$; these
functions remain bounded. The same choice of $N$ therefore controls
any prescribed number of parameter derivatives.

In the complementary region $|\eta|\geq\tau^{1/2}$, use the
weak resolvent estimates from the chapter on pseudodifferential
operators. For every prescribed set of derivatives, the remainders
are bounded by a fixed power of $1+\tau$ multiplied by a
negative power of $\langle\eta\rangle$ whose exponent can be made as large in magnitude as
desired by increasing the length of the parametrix. Conversion to a
joint bound is seen in the inequality
\[
(1+\tau)^a\langle\eta\rangle^{-J}
\leq C_{a,K}\langle(\tau,\eta)\rangle^{-K},
\qquad
|\eta|\geq\tau^{1/2},\quad \tau\geq1,\quad
J\geq2(a+K),
\]
where $a,K\geq0$. If $|\eta|\leq\tau$, use
$\langle\eta\rangle^{-J}\leq\tau^{-J/2}$; if $|\eta|\geq\tau$,
absorb $(1+\tau)^a$ into $\langle\eta\rangle^a$. The order of the
almost analytic extension is increased after fixing $J$, so that the
powers of $|t|^{-1}$ that have appeared can be integrated. The derivatives of
$Q_\tau$ lie in this same region, so the two localizations
can be recombined with the stated bounds. Finally, for
$R^{-1}\leq\tau\leq1$ the horizontal parameter stays bounded;
smoothing spectral calculus in the vertical variable and horizontal
localization give jointly smooth remainders. This justifies, for the
remainders as well, that in the extreme cones the symbol is respectively
$I$ and $0$ modulo arbitrarily negative order.

More precisely, after localization in $x$ and in a chart of $Z$,
we obtain a symbol $f$ with homogeneous components $f_{-j}$, extended
using a common cutoff near the origin. For $J\in\mathbb N$ and
multiindices $\alpha,\beta\in\mathbb N_0^{\dim Z}$ and
$\delta\in\mathbb N_0^n$ we have
\begin{equation}
\left\|\partial_z^\alpha\partial_\eta^\beta\partial_\xi^\delta
\left(f-\sum_{j=0}^{J-1}f_{-j}\right)\right\|_{\operatorname{op}}
\leq C_{J,\alpha,\beta,\delta,R}
\langle(\xi,\eta)\rangle^{-J-|\beta|-|\delta|}.
\label{eq:resto-simbolo-conjunto-corte}
\end{equation}
The same estimate allows derivatives in $x$ after inserting the
cutoffs. The smooth-kernel remainders satisfy these bounds for every
$J$: use the identity with $S_N$ between vertical Sobolev spaces
of arbitrary orders and apply the kernel regularity established
in Chapter~\ref{cap:operadores-pseudodiferenciales}. Derivatives in
$\xi$ are obtained by differentiating the resolvent and $\tau$ before integration.
The principal term, upon replacing $1+|\xi|^2$ by its homogeneous part,
is \eqref{eq:simbolo-principal-corte-espectral}. Its extensions to the
axes are smooth because $h$ is constant on the specified intervals.

In Kohn--Nirenberg quantization, $P$ acts only in $x$ and the symbol
of the cutoff is independent of $x$ before localization. The symbol of
$(P\otimes I)\Phi_R$ is therefore the product $p(x,\xi)f(z,\eta,\xi)$.
In the cone where $|\eta|$ dominates $|\xi|$, the second factor has arbitrarily
negative order; in its complement, $|\xi|$ controls the
joint frequency. The symbol estimates for $p$ and
\eqref{eq:resto-simbolo-conjunto-corte} prove that the product has
joint order one and is classical. The assertion for
$\Phi_R(P\otimes I)$ follows by taking adjoints.

For $D(I-\Phi_R)$ use vertical composition with a parameter.
In the cone where $|\xi|$ dominates $|\eta|$, the symbol of $I-\Phi_R$
has arbitrarily negative order. In the remaining cone, $|\eta|$
controls the joint frequency, and the symbol estimates
for $D$ apply. This yields a classical joint symbol of order one with principal
part $b_1(I-f_R)$. The localizations and changes of coordinates
are those of the classical calculus already constructed, completing the proof.
\end{proof}

\begin{lemma}[Product symbol and normal Thom class]
\label{lem:simbolo-producto-thom-normal}
Let $M$ be a closed smooth manifold and let
$a=[\mathbf{E}^+,\mathbf{E}^-,p]\in K_c^0(T^*M)$, where
$\mathbf{E}^+,\mathbf{E}^-\to M$ are smooth complex vector bundles.
Use the normal connection, metrics, and compactification
$\nu\subset\widehat\nu$ fixed above. Denote by
$q:T^*\nu\to M$ the projection and set
\[
\boldsymbol{\Lambda}_\nu^+
=q^*\!\left(\bigoplus_{\substack{j\in\{0,\ldots,r\}\\j\ \mathrm{even}}}
\Lambda^j(\nu\otimes_{\mathbb R}\mathbb C)\right),
\qquad
\boldsymbol{\Lambda}_\nu^-
=q^*\!\left(\bigoplus_{\substack{j\in\{0,\ldots,r\}\\j\ \mathrm{odd}}}
\Lambda^j(\nu\otimes_{\mathbb R}\mathbb C)\right).
\]
The component from even to odd degree of the normal Koszul homomorphism is
\[
\beta^+(v,\eta)
=\left.\bigl(\varepsilon(v+i\eta^\sharp)+
\iota(v-i\eta^\sharp)\bigr)\right|_{\boldsymbol{\Lambda}_\nu^+}.
\]
On $T^*\nu$, the graded matrix
\begin{equation}
\mathbf{c}(x,v;\xi,\eta)=
\begin{pmatrix}
 p(x,\xi)\otimes I&-I\otimes\beta^+(v,\eta)^*\\
 I\otimes\beta^+(v,\eta)&p(x,\xi)^*\otimes I
\end{pmatrix},
\label{eq:matriz-producto-thom-normal}
\end{equation}
is invertible outside a compact set. Its class, extended from this open subset,
satisfies
\begin{equation}
\operatorname{ext}_{T^*\nu}^{T^*\widehat\nu}[\mathbf{c}]
=
\operatorname{ext}_{T^*\nu}^{T^*\widehat\nu}
\operatorname{Th}_{T\nu}\bigl((\sharp_g)_*a\bigr)
\quad\text{in }K_c^0(T^*\widehat\nu).
\label{eq:clase-producto-thom-normal}
\end{equation}
\end{lemma}

\begin{proof}
The connection identifies the normal bundle $T\nu$ over $TM$ with
$\pi_{TM}^*(\nu\oplus\nu)$. The metric identifies the second copy of
$\nu$ with the vertical covariables, and the complex structure of
\eqref{eq:estructura-compleja-normal-tm} expresses its coordinate as
\[
z=v+i\eta^\sharp,
\qquad
J(v,\eta^\sharp)=(-\eta^\sharp,v).
\]
The normal Koszul homomorphism on the full exterior bundle is
$\beta(v,\eta)=\varepsilon(z)+\iota(\overline z)$. It is self-adjoint and
odd and satisfies
\[
\beta^2=(|v|^2+|\eta|^2)I,
\qquad
(\beta^+)^*\beta^+=(|v|^2+|\eta|^2)I
\quad\text{on }\boldsymbol{\Lambda}_\nu^+.
\]
The corresponding identity holds in odd degree. Thus $\beta^+$ is
invertible away from the zero section. These formulas are preserved under
changes of frame of $O(r)$ and define global homomorphisms of the stated
bundles.
The graded product of the horizontal complex $(\mathbf{E}^+,\mathbf{E}^-,p)$ with this
normal complex has even and odd parts, respectively,
\[
(q^*\mathbf{E}^+\otimes \boldsymbol{\Lambda}_\nu^+)
\oplus(q^*\mathbf{E}^-\otimes \boldsymbol{\Lambda}_\nu^-),
\qquad
(q^*\mathbf{E}^-\otimes \boldsymbol{\Lambda}_\nu^+)
\oplus(q^*\mathbf{E}^+\otimes \boldsymbol{\Lambda}_\nu^-),
\]
and its odd differential from even to odd degree is exactly the matrix
\eqref{eq:matriz-producto-thom-normal}. By the definition of the product of
complexes and the Thom isomorphism, its class is
$\operatorname{Th}_{T\nu}((\sharp_g)_*a)$, identifying $T\nu$ with
$T^*\nu$ by the fixed metrics. To check the support,
the product of the matrix with its adjoint has diagonal blocks
\[
p^*p\otimes I+(|v|^2+|\eta|^2)I,
\qquad
pp^*\otimes I+(|v|^2+|\eta|^2)I.
\]
Thus the matrix can fail to be invertible only when
$v=\eta=0$ and $p(x,\xi)$ is not invertible. The latter set is
compact because $a$ is a compactly supported class.

The factors are ordered with the horizontal factor first and the normal factor second; the sign
of the upper-right block is the corresponding Koszul sign and does not
introduce an additional sign into the class. In complex rank one, the
preceding coordinate gives the convention
$\lambda_{\mathbb C}=-b$ of
\eqref{eq:thom-linea-menos-bott}; by multiplicativity, it fixes both the
orientation and the sign in every rank. Finally, the matrix class has
compact support in $T^*\nu$. The extension-by-zero map for the
open subset $T^*\nu\subset T^*\widehat\nu$ gives, by definition, precisely
the right-hand side of \eqref{eq:clase-producto-thom-normal}.
\end{proof}

Spectral cutoffs allow the product to be formed within the
classical calculus on $\widehat\nu$. The vertical kernel provides a
copy of the sections over $M$. On its complement, a symmetry
interchanging the degrees will make the contributions to the index cancel.

\begin{lemma}[Horizontal and vertical product]
\label{lem:producto-horizontal-vertical}
Let $M$ be a closed smooth manifold and let
$\mathbf{E}^+,\mathbf{E}^-\to M$ be smooth complex vector bundles.
Let $P\in\Psi^1_{\mathrm{cl}}(M;\mathbf{E}^+,\mathbf{E}^-)$ be elliptic and
let $B$ be the operator of Lemma~\ref{lem:representante-normal-equivariante}.
There exist a classical elliptic operator of order one
\[
C:\Gamma(\widehat\nu,\boldsymbol{\mathcal E}^+)
\longrightarrow\Gamma(\widehat\nu,\boldsymbol{\mathcal E}^-)
\]
and bundles
\begin{align*}
\boldsymbol{\mathcal E}^+
&=(\widehat\pi^*\mathbf{E}^+\otimes\mathbf{W}_\nu^+)
\oplus(\widehat\pi^*\mathbf{E}^-\otimes\mathbf{W}_\nu^-),\\
\boldsymbol{\mathcal E}^-
&=(\widehat\pi^*\mathbf{E}^-\otimes\mathbf{W}_\nu^+)
\oplus(\widehat\pi^*\mathbf{E}^+\otimes\mathbf{W}_\nu^-)
\end{align*}
such that, if $a=[\boldsymbol\sigma(P)]$, then
\begin{equation}
[\boldsymbol\sigma(C)]
=\operatorname{ext}_{T^*\nu}^{T^*\widehat\nu}
\operatorname{Th}_{T\nu}\bigl((\sharp_{\mathbf{g}})_*a\bigr)
\quad\text{in }K_c^0(T^*\widehat\nu)
\label{eq:clase-cr-es-thom-normal}
\end{equation}
and
\begin{equation}
\operatorname{ind}_{\mathrm a}(C)=\operatorname{ind}_{\mathrm a}(P).
\label{eq:indice-producto-horizontal-vertical}
\end{equation}
The domain of the closure of $C$ is $H^1(\widehat\nu,\boldsymbol{\mathcal E}^+)$,
and that of its Hilbert space adjoint is
$H^1(\widehat\nu,\boldsymbol{\mathcal E}^-)$.
\end{lemma}

\begin{proof}
Set $\mathbf{W}=\mathbf{W}^+\oplus\mathbf{W}^-$ and consider on the
fiber the operator and grading
\[
D_B=\begin{pmatrix}0&B^*\\B&0\end{pmatrix},
\qquad
\Gamma_B=\begin{pmatrix}I&0\\0&-I\end{pmatrix}.
\]
By Lemma~\ref{lem:estimacion-conica-horizontal-vertical-parametro},
$D_B$ is self-adjoint with domain $H^1(S^r,\mathbf{W})$. Denote by
$\Pi_0$ the orthogonal projection onto its kernel. Elliptic regularity
shows that $\Pi_0$ is smoothing. By
\eqref{eq:nucleos-bott-compactificado}, its image is a line in even
degree, and the action of $G$ on it is trivial. Choose a generator
$\mathbf{e}_0$ with $L^2(S^r,\mathbf W)$ norm equal to one. Both $D_B$ and $\Pi_0$ are
$G$--equivariant and therefore induce vertical operators on
$\widehat\nu$; we retain the same letters for them.

Equip $\widehat\nu$ with the connection metric making
horizontal and vertical directions orthogonal, using the round metric on
the fibers. In an associated trivialization, its measure is
$d\lambda_{\mathbf{g}}(x)\,d\lambda_{S^r}(z)$. Changes of
trivialization are fiber isometries. Thus
\begin{equation}
(\mathcal I\mathbf{u})(x,z)=\mathbf{u}(x)\otimes\mathbf{e}_0(z)
\label{eq:inclusion-nucleo-vertical-base}
\end{equation}
defines a global isometry from $L^2(M,\mathbf{E}^+\oplus\mathbf{E}^-)$
onto $\Pi_0L^2(\widehat\nu,\boldsymbol{\mathcal E})$, where
$\boldsymbol{\mathcal E}=\boldsymbol{\mathcal E}^+\oplus
\boldsymbol{\mathcal E}^-$. It also identifies
$H^1(M,\mathbf E^+\oplus\mathbf E^-)$ with
$\Pi_0H^1(\widehat\nu,\boldsymbol{\mathcal E})$: in a product chart, derivatives in $x$ act on
$\mathbf{u}$ and vertical derivatives on the fixed smooth section
$\mathbf{e}_0$; the preceding isometry, also applied to base derivatives, gives the reverse bound.
Finitely many charts give the global equivalence.

First suppose that $n=\dim M\geq1$. Choose a finite cover
$(U_\ell)_{\ell=1}^J$ trivializing $\nu$ and $\mathbf{E}^\pm$, and
real functions $\chi_\ell\in C_c^\infty(U_\ell)$ such that
\[
\sum_{\ell=1}^J\chi_\ell^2=1\quad\text{on }M.
\]
Use orthonormal frames for $\nu$ and unitary frames for the complex
bundles. When passing to Euclidean coordinates in $U_\ell$, incorporate the
square root of the measure Jacobian into the unitary isomorphisms
from square-integrable sections of $\mathbf E^+\restriction_{U_\ell}$ and
$\mathbf E^-\restriction_{U_\ell}$ to their coordinate representations on
$\phi_\ell(U_\ell)\subseteq\mathbb R^n$, where $\phi_\ell$ is the chosen chart. Thus the adjoints in coordinates are the
Hilbert space adjoints of the transported operators.

Let $P_\ell$ be the local expression of $\zeta_\ell P\zeta_\ell$, extended
to $\mathbb R^n$, where $\zeta_\ell\in C_c^\infty(U_\ell)$ equals
one near $\operatorname{supp}\chi_\ell$. Set
\[
A_\ell=\begin{pmatrix}0&P_\ell^*\\P_\ell&0\end{pmatrix},
\qquad
\Gamma_P=\begin{pmatrix}I&0\\0&-I\end{pmatrix}
\]
on the factor $\mathbf{E}^+\oplus\mathbf{E}^-$. In the product chart
$U_\ell\times S^r$ apply
Lemma~\ref{lem:cortes-espectrales-producto-clasico} to $D_B$ and write
$\Phi_{\ell,R}$ for the cutoff constructed there. Before multiplying by
$\chi_\ell$, perform horizontal quantization on all of $\mathbb R^n$.
Define the odd operator
\begin{equation}
\mathscr C_R
=\sum_{\ell=1}^J\chi_\ell
\left[
\frac12\bigl(A_\ell\Phi_{\ell,R}+\Phi_{\ell,R}A_\ell\bigr)
+R\Gamma_P D_B(I-\Phi_{\ell,R})
\right]\chi_\ell.
\label{eq:producto-operador-normal-bott}
\end{equation}
Transport each summand to $\widehat\pi^{-1}(U_\ell)$ and extend it
by zero in the base variable. The outer functions make this
extension globally defined. The cutoff lemma proves that
$\mathscr C_R\in\Psi^1_{\mathrm{cl}}(\widehat\nu;
\boldsymbol{\mathcal E},\boldsymbol{\mathcal E})$. Symmetrization of the
horizontal term and commutation of $D_B$ with the cutoffs show
that it is formally self-adjoint.

Let us verify its ellipticity and fix the size of $R$. Denote by
$a_1(x,\xi)$ the odd principal matrix of $P$ and by $b_1(y,\eta)$ that
of $D_B$, where $y\in\widehat\nu$, $x=\widehat\pi(y)$, and
$(\xi,\eta)\in\mathbf{H}_y^*\oplus\mathbf{V}_y^*$. Extend both
symbols by zero at their respective origins. They are continuous and homogeneous of
degree one. There exist constants $c_P,c_B,C_B>0$ such that
\begin{equation}
a_1(x,\xi)^2\geq c_P^2|\xi|^2I,
\qquad
c_B^2|\eta|^2I\leq b_1(y,\eta)^2\leq C_B^2|\eta|^2I.
\label{eq:positividad-simbolo-producto-horizontal-vertical}
\end{equation}
Homogeneity also gives a constant $L_P>0$ for which
\[
\|a_1(x,\xi')-a_1(x,\xi)\|_{\operatorname{op}}
\leq L_P|\xi'-\xi|.
\]
Indeed, first derivatives in the covariable are homogeneous of
degree zero and bounded in a finite atlas. Integrate this bound along
the segment between $\xi$ and $\xi'$; if it passes through the origin,
apply it on the two subsegments and use continuity there.

The product coordinates of the trivialization and the fixed horizontal
splitting need not give the same decomposition. In the
trivialization $\ell$, the coordinate horizontal covariable corresponds to
\[
\xi_\ell=\xi+T_\ell(y)\eta\in T_x^*M,
\qquad \|T_\ell(y)\|_{\operatorname{op}}\leq A_0.
\]
The constant $A_0$ exists because we use only finitely many
trivializations over compact sets. If $|\cdot|_\ell$ is the Euclidean
norm of the components in that chart, fix $c_x,C_x>0$ such that
\[
c_x|\omega|\leq|\omega|_\ell\leq C_x|\omega|,
\qquad
\omega\in T_x^*M,\quad x\in\operatorname{supp}\chi_\ell,
\quad\ell\in\{1,\ldots,J\}.
\]
The symbol of the cutoff, transported to the corresponding fiber, is
\[
f_{\ell,R}(y;\xi,\eta)
=h\!\left(\frac{R^2b_1(y,\eta)^2}{|\xi_\ell|_\ell^2}\right).
\]
On the axes, use the extensions of the lemma. These endomorphisms
are positive, bounded by $I$, and commute with one another and with $b_1$,
since they are functions of the same matrix $b_1^2$. Set
\[
F_R=\sum_{\ell=1}^J\chi_\ell(x)^2f_{\ell,R},
\]
extending each term by zero outside its trivialization.
The principal symbol of \eqref{eq:producto-operador-normal-bott} is
\begin{equation}
c_R
=\sum_{\ell=1}^J\chi_\ell^2
a_1(x,\xi_\ell)f_{\ell,R}
+R\Gamma_Pb_1(y,\eta)(I-F_R).
\label{eq:simbolo-producto-horizontal-vertical-regularizado}
\end{equation}
The identity factors on $\mathbf{E}$ and $\mathbf{W}$ remain
implicit in this formula; the two operators act on distinct tensor
factors.

Compare $c_R$ with
\[
m_R=a_1(x,\xi)F_R+R\Gamma_Pb_1(y,\eta)(I-F_R).
\]
The preceding bounds imply
\begin{equation}
\|c_R-m_R\|_{\operatorname{op}}
\leq L_PA_0|\eta|.
\label{eq:error-coordenadas-producto-normal}
\end{equation}
Since $a_1\Gamma_P=-\Gamma_Pa_1$, and $F_R$ commutes with the vertical
factors, the cross terms in the square cancel exactly:
\begin{equation}
m_R^2=a_1(x,\xi)^2F_R^2
+R^2b_1(y,\eta)^2(I-F_R)^2.
\label{eq:estimacion-producto-normal}
\end{equation}

Choose
\[
\delta=\min\left\{\frac14,\frac{c_x}{4C_B},\frac{c_B}{8C_x}\right\},
\qquad R\geq\max\{1,A_0/\delta\}.
\]
If $R|\eta|\leq\delta|\xi|$, then
$|\xi_\ell|\geq(1-\delta^2)|\xi|\geq\frac34|\xi|$.
All matrices $R^2b_1^2/|\xi_\ell|_\ell^2$ have spectrum in
$[0,1]$ for the indices
$\ell\in\{1,\ldots,J\}$ with $\chi_\ell(x)\neq0$. Thus
$f_{\ell,R}=I$ in the terms occurring in the sum, and $F_R=I$.
If $|\xi|\leq\delta R|\eta|$, then
$|\xi_\ell|_\ell\leq2C_x\delta R|\eta|$.
The corresponding spectra then lie in $[4,+\infty)$;
hence $F_R=0$. In the remaining region, we have simultaneously
\[
|\xi|>\delta R|\eta|,
\qquad R|\eta|>\delta|\xi|.
\]
Thus each of $|\xi|^2$ and $R^2|\eta|^2$ is at least
$\frac{\delta^2}{2}(|\xi|^2+R^2|\eta|^2)$. Moreover,
\[
F_R^2+(I-F_R)^2
=2(F_R-\tfrac12I)^2+\tfrac12I\geq\tfrac12I.
\]
Applying \eqref{eq:estimacion-producto-normal} in the three regions
gives, with $\displaystyle c_*=\min\{c_P,c_B\}$,
\[
\|m_R\mathbf{v}\|\geq
\frac{\delta c_*}{2}
(|\xi|^2+R^2|\eta|^2)^{1/2}\|\mathbf{v}\|
\]
for every fiber vector $\mathbf{v}$. Finally, fix
\[
R\geq
\max\left\{1,\frac{A_0}{\delta},
\frac{4L_PA_0}{\delta c_*}\right\}.
\]
Inequality \eqref{eq:error-coordenadas-producto-normal} implies
\begin{equation}
\|c_R\mathbf{v}\|\geq
\frac{\delta c_*}{4}
(|\xi|^2+R^2|\eta|^2)^{1/2}\|\mathbf{v}\|.
\label{eq:positividad-uniforme-matriz-conica}
\end{equation}
This proves ellipticity of $\mathscr C_R$.
Lemma~\ref{lem:estimacion-conica-horizontal-vertical-parametro} and
Theorem~\ref{pseudo:teo-fredholm-eliptico} give
\begin{equation}
\mathcal D(\mathscr C_R)=H^1(\widehat\nu,\boldsymbol{\mathcal E}),
\qquad
\mathcal D(\mathscr C_R^\pm)=H^1(\widehat\nu,\boldsymbol{\mathcal E}^\pm).
\label{eq:dominio-producto-normal}
\end{equation}
The operator $\mathscr C_R$ is self-adjoint, and its component
$C:=\mathscr C_R^+$ is Fredholm, with adjoint $\mathscr C_R^-$ on the
stated domain.

Let us calculate its index. The operator $D_B^2+\Pi_0$ is positive, invertible,
and classical elliptic of order two on the fiber. The complex powers of
Theorem~\ref{pseudo:teo-potencias-complejas-clasicas} allow us to define
\[
J_B=D_B(D_B^2+\Pi_0)^{-1/2}\in\Psi^0_{\mathrm{cl}}(S^r;\mathbf{W},\mathbf{W}).
\]
It is equivariant and satisfies
\[
J_B^*=J_B,\qquad J_B^2=I-\Pi_0,
\qquad J_B\Gamma_B=-\Gamma_BJ_B.
\]
The vertical operators $\Pi_0$ and $J_B$ are bounded on the spaces
$L^2(\widehat\nu,\boldsymbol{\mathcal E})$ and $H^1(\widehat\nu,\boldsymbol{\mathcal E})$ of the total manifold. To check the second assertion,
in a product chart they commute with derivatives in $x$ and are bounded on
vertical $H^1(S^r,\mathbf W)$ by the order-zero calculus; a finite partition of
unity gives the global bound.

Every summand of \eqref{eq:producto-operador-normal-bott} commutes with
$\Pi_0$. On the complement of its image define
\[
\mathcal T=\Gamma_P\Gamma_BJ_B.
\]
This operator is unitary, interchanges degrees, and satisfies
$\mathcal T^2=-I$. The cutoffs commute with $J_B$ and $\Gamma_B$;
the operators $A_\ell$ commute with them and anticommute with
$\Gamma_P$. Thus
\begin{equation}
\mathcal T\mathscr C_R=-\mathscr C_R\mathcal T
\quad\text{on }(I-\Pi_0)H^1(\widehat\nu,\boldsymbol{\mathcal E}).
\label{eq:simetria-complemento-nucleo-vertical}
\end{equation}
First verify the identity on smooth sections. Boundedness of
$\mathcal T$ on $H^1(\widehat\nu,\boldsymbol{\mathcal E})$ and \eqref{eq:dominio-producto-normal} allows
it to extend to the domain. Consequently, $\mathcal T$ identifies the
even and odd parts of the kernel of $\mathscr C_R$ lying in
$(I-\Pi_0)L^2(\widehat\nu,\boldsymbol{\mathcal E})$. Their dimensions are finite and their contributions to the
index cancel.

On $\Pi_0L^2(\widehat\nu,\boldsymbol{\mathcal E})$, the vertical operator vanishes and the cutoff is the
scalar horizontal multiplier
\[
a_{0,R}(\xi)=h\!\left(\frac{R^2}{1+|\xi|^2}\right).
\]
This function equals one outside a ball. Through
\eqref{eq:inclusion-nucleo-vertical-base}, the even-to-odd component of
the restriction of $\mathscr C_R$ is identified with the operator
\[
P_R=\frac12\sum_{\ell=1}^J\chi_\ell
\bigl(P_\ell a_{0,R}(D_x)+a_{0,R}(D_x)P_\ell\bigr)\chi_\ell
\]
on $M$, transporting each summand from its chart. Since
$a_{0,R}-1$ has compact support in the horizontal frequency,
\begin{equation}
P_R-P=\sum_{\ell=1}^J\chi_\ell[P,\chi_\ell]+S_R,
\qquad S_R\in\Psi^{-\infty}(M;\mathbf{E}^+,\mathbf{E}^-).
\label{eq:restriccion-producto-nucleo-vertical}
\end{equation}
The commutators are classical of order zero. By
Lemma~\ref{lem:compacidad-orden-cero-h1-l2}, $P_R-P$ is compact from
$H^1(M,\mathbf E^+)$ to $L^2(M,\mathbf E^-)$, so
$\operatorname{ind}P_R=\operatorname{ind}P$ by
Theorem~\ref{teo:estabilidad-indice-fredholm}.

If $K^\pm$ denotes the part of the kernel of $\mathscr C_R$ lying in
$(I-\Pi_0)L^2(\widehat\nu,\boldsymbol{\mathcal E})$ and of degree $\pm$, we have proved the decompositions
\begin{equation}
\begin{aligned}
\ker C&=\mathcal I(\ker P_R)\oplus K^+,\\
\ker C^*&=\mathcal I(\ker P_R^*)\oplus K^-,
\qquad \mathcal T:K^+\xrightarrow{\;\cong\;}K^-.
\end{aligned}
\label{eq:nucleos-producto-normal}
\end{equation}
Subtracting dimensions yields
\begin{equation}
\operatorname{ind}C=\operatorname{ind}P_R=\operatorname{ind}P.
\label{eq:indice-homotopia-conica}
\end{equation}

It remains to identify the symbol class. In
\eqref{eq:simbolo-producto-horizontal-vertical-regularizado} replace
$T_\ell$ simultaneously by $tT_\ell$, with $t\in[0,1]$, both in
$a_1(x,\xi_\ell)$ and in $f_{\ell,R}$. The estimates leading to
\eqref{eq:positividad-uniforme-matriz-conica} remain valid,
since $\|tT_\ell\|_{\operatorname{op}}\leq A_0$. This gives a
homotopy of invertible symbols on the cotangent bundle with its zero section removed.
At $t=0$ it takes the form
\[
a_1(x,\xi)F+R\Gamma_Pb_1(y,\eta)(I-F),
\]
where $F$ is a convex sum of functions of $b_1^2$ and $0\leq F\leq I$.
For $s\in[0,1]$ replace the two coefficients by
\[
F_s=(1-s)F+sI,\qquad G_s=(1-s)(I-F)+sI.
\]
The square of $a_1F_s+R\Gamma_Pb_1G_s$ is
$a_1^2F_s^2+R^2b_1^2G_s^2$. Since $F_s\geq F$ and $G_s\geq I-F$,
the lower bound at the start of the homotopy is preserved. At the end we
obtain $a_1+R\Gamma_Pb_1$. Deforming the positive factor $R$ to
one yields $a_1+\Gamma_Pb_1$.

These last homotopies are continuous homotopies of bundle homomorphisms on the
cosphere bundle. They are not used to define operators at intermediate
values. Continuity suffices for homotopy invariance of the
$K$--theory class. The even-to-odd component of the final symbol is
\[
\begin{pmatrix}
p_1(x,\xi)\otimes I&-I\otimes b_1^+(y,\eta)^*\\
I\otimes b_1^+(y,\eta)&p_1(x,\xi)^*\otimes I
\end{pmatrix},
\]
where $b_1^+$ is the symbol of $B$. This is the graded product of the horizontal
class with the vertical class.
Lemma~\ref{lem:representante-normal-equivariante} identifies the latter
with the extension of the normal Bott class.
Lemma~\ref{lem:simbolo-producto-thom-normal} identifies its product with
the Thom class, with the order of factors and orientation already fixed.
This gives \eqref{eq:clase-cr-es-thom-normal}.

If $\dim M=0$, $M$ is a finite set. Over each point, use
the vertical operator $\Gamma_PD_B$ directly in the graded product
with $\mathbf{E}^+\oplus\mathbf{E}^-$. Its index is
$\dim\mathbf{E}^+-\dim\mathbf{E}^-$, which agrees with the index of
any linear map from $\mathbf{E}^+$ to $\mathbf{E}^-$ over that
point. Its symbol is the product of this virtual class with the Bott
class. Summing over the points of $M$ gives both identities in the statement
and completes the proof.
\end{proof}

\begin{theorem}[Analytic compatibility with the pushforward]
\label{teo:compatibilidad-indice-pushforward}
Let $M$ be a closed smooth manifold, let
$j\colon M\hookrightarrow\mathbb R^N$ be a smooth embedding, and let
$a\in K_c^0(T^*M)$. Then
\begin{equation}
\operatorname{ind}_{\mathrm a}^{\mathbb R^N}
(j_!a)
=\operatorname{ind}_{\mathrm a}^{M}(a).
\label{eq:compatibilidad-indice-pushforward}
\end{equation}
\end{theorem}

\begin{proof}
Represent $a$ by the symbol of an elliptic operator $P$ on $M$.
Theorem~\ref{teo:realizacion-clases-simbolo} allows this after
stabilization, which adds an identity of index zero. If
the resulting representative has order zero, fix a positive invertible
operator
\[
\Lambda_+\in\Psi^1_{\mathrm{cl}}(M;\mathbf{E}^+),
\qquad
\Lambda_+^{-1}\in\Psi^{-1}_{\mathrm{cl}}(M;\mathbf{E}^+),
\]
and replace $P$ by $P\Lambda_+$. Composition with $\Lambda_+$ does not
change the index. Its positive scalar symbol contracts to the identity
symbol after normalizing the order, so it also leaves the class
$a$ unchanged. We may therefore apply
Lemma~\ref{lem:producto-horizontal-vertical}. This gives an operator $C$ on
$\widehat\nu$ with
\begin{equation}
\operatorname{ind}_{\mathrm a}(C)=\operatorname{ind}_{\mathrm a}(P)
\label{eq:indice-c-sobre-normal-compactificado}
\end{equation}
and symbol class
\begin{equation}
[\boldsymbol{\sigma}(C)]
=\operatorname{ext}_{T^*\nu}^{T^*\widehat\nu}
\operatorname{Th}_{T\nu}\bigl((\sharp_g)_*a\bigr)
\in K_c^0(T^*\widehat\nu),
\label{eq:simbolo-c-es-thom}
\end{equation}
where the metric and connection identify the Thom representative on
$T\nu$ with the corresponding complex on $T^*\nu$, followed by extension
from the open subset $T^*\nu\subset T^*\widehat\nu$. This equality is
precisely \eqref{eq:clase-cr-es-thom-normal}.

To pass from $\widehat\nu$ to the tubular neighborhood, we reduce
the order of the operator and localize its symbol.
Let
\[
\Lambda_{-1}\in\Psi^{-1}_{\mathrm{cl}}(\widehat\nu;\boldsymbol{\mathcal E}^+)
\]
be positive and elliptic, chosen so that
$\Lambda_{-1}:L^2(\widehat\nu,\boldsymbol{\mathcal E}^+)
\longrightarrow H^1(\widehat\nu,\boldsymbol{\mathcal E}^+)$ is an isomorphism; its
inverse has order one. The operator
\[
 C\Lambda_{-1}\colon L^2(\widehat\nu,\boldsymbol{\mathcal E}^+)
 \longrightarrow L^2(\widehat\nu,\boldsymbol{\mathcal E}^-)
\]
has order
zero. Since $\Lambda_{-1}$ is invertible, it preserves the kernel, cokernel, and
index of $C$; its normalized positive symbol is homotopic to the
identity and also leaves the symbol class unchanged. Denote by $A$ an
order-zero quantization obtained in this way. By
\eqref{eq:simbolo-c-es-thom}, its symbol triple is
elementary on a neighborhood of $\mathbf{s}_\infty$. After adding an elementary
triple, the input and output bundles are identified on that neighborhood.
The polar retraction
\[
s\longmapsto s(s^*s)^{-\frac{1}{2}}
\]
is applied only to the symbol isomorphism in that region and makes it
unitary. The triple being elementary means precisely that, after the
chosen stabilization, this unitary admits a homotopy there to the identity
through isomorphisms. With a single cutoff supported in the same region,
quantization of this homotopy is continuous in
\[
 \mathcal L\bigl(
 L^2(\widehat\nu,\boldsymbol{\mathcal E}^+),
 L^2(\widehat\nu,\boldsymbol{\mathcal E}^-)\bigr);
\]
and all its members are elliptic and
Fredholm, so the index remains constant.

At the end of the homotopy, replacing $A$ by the final quantization, the
principal symbol is the identity on a neighborhood $\mathcal U_0$ of
$\mathbf{s}_\infty$. Fix a smaller neighborhood $\mathcal U_\infty$ whose closure
is contained in $\mathcal U_0$. The identification already chosen between
$\boldsymbol{\mathcal E}^+$ and $\boldsymbol{\mathcal E}^-$ over $\mathcal U_0$ determines an
identity morphism $I_\infty$ there. Choose
$\zeta\in C_c^\infty(\mathcal U_0)$ equal to one on a neighborhood of
$\overline{\mathcal U_\infty}$ and extend the morphism
$\zeta I_\infty$ by zero. This gives a global smooth morphism
\[
 R:\boldsymbol{\mathcal E}^+\longrightarrow\boldsymbol{\mathcal E}^-
\]
agreeing with $I_\infty$ near
$\overline{\mathcal U_\infty}$; we do not assert that $R$ is invertible outside
this region. If $r=\boldsymbol{\sigma}_0^\Psi(R)$, let
\[
K_A:=\operatorname{supp}\bigl(\boldsymbol{\sigma}_0^\Psi(A)-r\bigr)
\subset S^*\widehat\nu.
\]
This actual support is compact, and its projection to $\widehat\nu$ does not meet
$\mathcal U_\infty$. Choose
$\chi\in C^\infty(\widehat\nu)$ equal to one on a neighborhood of
$\pi(K_A)$ and zero on $\mathcal U_\infty$, and define
\begin{equation}
A_0:=R+\chi(A-R)\chi.
\label{eq:localizacion-operador-normal}
\end{equation}
The operators $A$ and $R$ act between the same pair
$\boldsymbol{\mathcal E}^+,\boldsymbol{\mathcal E}^-$, and $\chi$ acts by scalar multiplication;
therefore, $A_0$
is a global operator on this bundle. Moreover,
\[
 \boldsymbol{\sigma}_0^\Psi(A_0)
 =r+\chi^2\bigl(\boldsymbol{\sigma}_0^\Psi(A)-r\bigr)
 =\boldsymbol{\sigma}_0^\Psi(A),
\]
because $\chi=1$ on the support of the difference. By composition
calculus, $A_0-A$ has order $-1$ on the closed manifold
$\widehat\nu$ and is compact as an operator
\[
 L^2(\widehat\nu,\boldsymbol{\mathcal E}^+)
 \longrightarrow L^2(\widehat\nu,\boldsymbol{\mathcal E}^-).
\]
On $\mathcal U_\infty$ we have
$\chi=0$ and $R=I_\infty$. Thus
\begin{equation}
\operatorname{ind}_{\mathrm a}(A_0)
=\operatorname{ind}_{\mathrm a}(A)
=\operatorname{ind}_{\mathrm a}(P),
\qquad
A_0=I_\infty\text{ on }\mathcal U_\infty.
\label{eq:indice-operador-normal-localizado}
\end{equation}
This formula supplies the exact equality near the boundary needed
to extend the operator; it is not the unbounded oscillator that is being extended.

After shrinking the tubular domain, we may assume that
$V\subset\nu$ is the open disk $|v|_{h_\nu}<\rho(x)$ for a positive smooth
function $\rho$. The proper radial contraction
\[
\kappa\colon\nu\longrightarrow V,
\qquad
\kappa(x,v)=\left(x,
\frac{\rho(x)v}{\sqrt{1+|v|_{h_\nu}^2}}\right),
\]
with a positive function $\rho$ chosen within the tubular radius sends the
regions $|v|_{h_\nu}\geq R$ to progressively thinner collars of the boundary of $V$.
Conjugate $A_0$ by the
unitary isomorphism incorporating the Jacobian of $\kappa$ and then by
the tubular diffeomorphism $\Theta\colon V\longrightarrow U$. Since $A_0$ is
exactly the identity near $\mathbf{s}_\infty$, the transported operator is
the identity on a collar of the boundary of $U$, under the identification already
chosen between the input and output bundles. Denote by
$d\in K_c^0(T^*U)$ the class of this localized symbol.

Extension must be performed after stabilizing the bundles; identification
on the collar alone does not suffice. Consider
\[
e:=\operatorname{ext}(d)\in K_c^0(T^*\mathbb R^N).
\]
In the difference model, choose finite-rank Hermitian complements of
the two bundles of the triple representing $e$ over a pair
$(B^{2N},S^{2N-1})$ whose interior contains the support. The space
$B^{2N}$ is contractible, so these complements and both bundles are
trivial. The isomorphism over $S^{2N-1}$ forces the input and output
bundles to have equal ranks. After adding the same trivial elementary
triple of sufficiently large rank, both are identified with
$\underline{\mathbb C}^{\,q}$. Stabilization adds the symbol of an
identity, so it changes neither the $K$--theory class nor the index.

On the collar, the identification coming from $A_0=I$ and the
preceding trivializations make the input and output bundle the same. The
polar retraction, followed by the elementary homotopy of the stabilized triple,
allows us to choose a representative
\[
s_e:T^*\mathbb R^N\longrightarrow
\operatorname{End}(\mathbb C^q)
\]
invertible outside a compact set in phase space, equal to $I$ when the
base variable lies outside a compact set, and whose restriction to $T^*U$
represents $d$. Quantize this homotopy with a single cutoff and localize
as in \eqref{eq:localizacion-operador-normal}; this gives
\[
\widetilde A\in\operatorname{Ell}_c^0
(\mathbb R^N;\mathbb C^q),
\qquad
[\boldsymbol{\sigma}(\widetilde A)]=e.
\]
All homotopies are norm-continuous on $L^2(\mathbb R^N,\mathbb C^q)$, and stabilization only
adds identities, so they do not alter the index. Now $A_0$, viewed on
$\widehat\nu$, and $\widetilde A$, viewed on $\mathbb R^N$, are the two
extensions of the same localized class $d$.
Proposition~\ref{prop:excision-indice-analitico} and
\eqref{eq:indice-operador-normal-localizado} give
\begin{equation}
\operatorname{ind}_{\mathrm a}^{\mathbb R^N}(\widetilde A)
=\operatorname{ind}_{\mathrm a}^{M}(P).
\label{eq:indice-extension-tubular}
\end{equation}
$e$ is the extension by zero of \eqref{eq:simbolo-c-es-thom}; by
\eqref{eq:pushforward-encaje-indice}, it is precisely $j_!a$. Since the
analytic index depends only on the symbol class,
\eqref{eq:indice-extension-tubular} proves
\eqref{eq:compatibilidad-indice-pushforward}.
\end{proof}

The preceding compatibility is the property of the analytic index that allows
passage from a closed manifold to its tubular neighborhood. The embedding
construction and its relation to the product of symbols are also developed in
\cite[equations~(12.8)--(12.16), pp.~293--301]{BleeckerBoossIndex}.

\section{The Atiyah--Singer theorem}

Before the main theorem, let us clarify why the analytic index can
be evaluated directly on a $K$--theory class.
Theorem~\ref{teo:realizacion-clases-simbolo} realizes every
$a\in K_c^0(T^*M)$ as the symbol of an elliptic operator, and
Theorem~\ref{teo:indice-analitico-sobre-k} proves that two realizations of
the same class have the same index. Indeed, after adding
identities, they are joined by a homotopy of elliptic symbols; this is
quantized as a continuous family of Fredholm operators, and
Theorem~\ref{pseudo:teo-estabilidad-indice} keeps the index constant.
This gives the homomorphism
\[
\operatorname{ind}_{\mathrm a}^M\colon K_c^0(T^*M)\longrightarrow\mathbb Z.
\]

\begin{theorem}[Atiyah--Singer, $K$--theoretic version]
\label{teo:atiyah-singer-k-cerradas}
\index{Atiyah--Singer theorem}
Let $M$ be a closed smooth manifold and let
$P\colon\Gamma(M,\mathbf{E}^+)\longrightarrow\Gamma(M,\mathbf{E}^-)$ be a classical elliptic
pseudodifferential operator between complex vector bundles. Then
\begin{equation}
\operatorname{ind}_{\mathrm a}(P)
=\operatorname{ind}_{\mathrm t}([\boldsymbol{\sigma}(P)])
\label{eq:atiyah-singer-k-cerradas}
\end{equation}
with the topological index of
Definition~\ref{def:indice-topologico-cerradas}.
\end{theorem}

\begin{proof}
Write $a=[\boldsymbol{\sigma}(P)]$. By
Theorem~\ref{teo:compatibilidad-indice-pushforward},
\begin{equation}
\operatorname{ind}_{\mathrm a}^{M}(a)
=\operatorname{ind}_{\mathrm a}^{\mathbb R^N}
(j_!a).
\label{eq:reduccion-indice-a-euclidiano}
\end{equation}
Theorem~\ref{teo:indice-euclidiano-k} and
\eqref{eq:normalizacion-pushforward-punto} state that, for every
$c\in K_c^0(T^*\mathbb R^N)$,
\[
\operatorname{ind}_{\mathrm a}^{\mathbb R^N}(c)
=(i_N)_!^{-1}(c).
\]
Applying this to $c=j_!a$ turns equality
\eqref{eq:reduccion-indice-a-euclidiano} into precisely
\[
\operatorname{ind}_{\mathrm a}^{M}(a)
=(i_N)_!^{-1}(j_!a)
=\operatorname{ind}_{\mathrm t}(a).
\]
\end{proof}

The argument also proves two properties often included in
the statement of the theorem. If $P_0$ and $P_1$ have symbols homotopic
through elliptic symbols, their indices agree. If a symbol
class is supported in $T^*U$ for an open subset $U\subseteq M$, its index
can be calculated in any manifold containing $U$. Neither
assertion requires a local formula for the heat kernel.

\section{Cohomological form}

The $K$--theoretic version is integral and does not require an orientation of $M$. To obtain an
integral formula, apply the compactly supported Chern character
\[
\operatorname{ch}_c\colon K_c^0(TM)\longrightarrow
H_c^{\mathrm{even}}(TM;\mathbb Q).
\]
The total space $TM$ has a canonical orientation, even if $M$ is
nonorientable: a connection identifies
\[
T(TM)\cong\pi_{TM}^*TM\oplus\pi_{TM}^*TM,
\]
and $J(X,Y)=(-Y,X)$ defines an almost complex structure. Denote by
$[TM]_{\mathrm{BM}}\in H_{2n}^{\mathrm{BM}}(TM;\mathbb Z)$ its Borel--Moore
fundamental class.

For a complex bundle $\mathbf{W}$ with formal Chern roots
$x_1,\ldots,x_s$, fix the normalizations
\begin{equation}
\operatorname{ch}(\mathbf{W})=\sum_{j=1}^s e^{x_j},
\qquad
\operatorname{Td}(\mathbf{W})=
\prod_{j=1}^s\frac{x_j}{1-e^{-x_j}}.
\label{eq:normalizaciones-chern-todd}
\end{equation}
In particular,
\[
\operatorname{Td}(\mathbf{W})
=1+\frac{1}{2}c_1(\mathbf{W})
+\frac{1}{12}\bigl(c_1(\mathbf{W})^2+c_2(\mathbf{W})\bigr)+\cdots.
\]

\begin{theorem}[Cohomological index formula]
\label{teo:formula-cohomologica-indice-cerradas}
Let $M$ be a closed smooth manifold of real dimension $n$,
let $\mathbf{E}^+,\mathbf{E}^-\to M$ be smooth complex vector bundles,
and let
$P\colon\Gamma(M,\mathbf{E}^+)\to\Gamma(M,\mathbf{E}^-)$ be a classical elliptic
pseudodifferential operator. Identify $T^*M$ with $TM$
using a metric and write
$a=(\sharp_g)_*[\boldsymbol{\sigma}(P)]$. Then
\begin{equation}
\operatorname{ind}_{\mathrm a}(P)
=(-1)^n
\left\langle
\operatorname{ch}_c(a)\smile
\pi_{TM}^*\operatorname{Td}(TM\otimes_{\mathbb R}\mathbb C),
[TM]_{\mathrm{BM}}
\right\rangle.
\label{eq:formula-cohomologica-indice-tm}
\end{equation}
If $M$ is oriented and
\[
\Phi_{TM}\colon H^*(M;\mathbb Q)\longrightarrow
H_c^{*+n}(TM;\mathbb Q)
\]
is the cohomological Thom isomorphism for this orientation, then
\begin{equation}
\operatorname{ind}_{\mathrm a}(P)
=(-1)^{\frac{n(n+1)}{2}}
\int_M
\left[
\Phi_{TM}^{-1}\!\left(\operatorname{ch}_c(a)\right)
\smile\operatorname{Td}(TM\otimes_{\mathbb R}\mathbb C)
\right]_{n}.
\label{eq:formula-cohomologica-indice-base}
\end{equation}
Here $[\cdot]_n$ denotes the component of degree $n$.
\end{theorem}

Compatibility between Thom and the Chern character is also proved
in \cite[Theorem~13.1, pp.~313--315]{BleeckerBoossIndex}.

\begin{proof}
Comparison of the Thom isomorphisms in $K$--theory and
cohomology is a purely topological result. With the normalizations
\eqref{eq:normalizaciones-chern-todd}, if $\mathbf{W}\longrightarrow X$ is a complex
bundle of rank $s$, $\lambda_{\mathbf{W}}\in K_c^0(\mathbf{W})$ is its Thom class, and
$\Phi_{\mathbf{W}}$ is the cohomological Thom isomorphism, then
\begin{equation}
\Phi_{\mathbf{W}}^{-1}\!\left(\operatorname{ch}_c(\lambda_{\mathbf{W}})\right)
=\prod_{j=1}^s\frac{1-e^{x_j}}{x_j}
=(-1)^s\operatorname{Td}(\overline{\mathbf{W}})^{-1}.
\label{eq:defecto-chern-thom}
\end{equation}
Appendix~\ref{ap:topologia} records this comparison with its reference;
the proof using the splitting principle and naturality of the Chern
character is in \cite[Theorem~13.1, pp.~313--315]{BleeckerBoossIndex}.

For the embedding $j\colon M\hookrightarrow\mathbb R^N$, the stable complex
sum
\[
\pi_{TM}^*(TM\otimes_{\mathbb R}\mathbb C)\oplus T\nu
\cong TM\times\mathbb C^N
\]
is trivial. By multiplicativity of the Todd class,
\begin{equation}
\operatorname{Td}(T\nu)
=\pi_{TM}^*\operatorname{Td}
(TM\otimes_{\mathbb R}\mathbb C)^{-1}.
\label{eq:todd-normal-inversa}
\end{equation}
Here both $T\nu$, with the structure of
\eqref{eq:estructura-compleja-normal-tm}, and
$TM\otimes_{\mathbb R}\mathbb C$ are isomorphic to their conjugates; hence
\eqref{eq:defecto-chern-thom} produces the Todd class appearing in
\eqref{eq:todd-normal-inversa}.
Since $T\nu$ has complex rank $r=N-n$, the Thom formula gives explicitly
\begin{equation}
\operatorname{ch}_c(j_!a)
=(-1)^r\operatorname{ext}\,\Phi_{T\nu}\!\left(
\operatorname{ch}_c(a)\smile
\pi_{TM}^*\operatorname{Td}
(TM\otimes_{\mathbb R}\mathbb C)
\right).
\label{eq:chern-pushforward-con-signo-normal}
\end{equation}
On the other hand, $\alpha_N(b^{\boxtimes N})=1$ and the normalization
$b=[\mathbf{E}_{-1}]-[1]$ satisfies
\[
\left\langle\operatorname{ch}_c(b^{\boxtimes N}),
[\mathbb C^N]_{\mathrm{BM}}\right\rangle=1.
\]
Thus \eqref{eq:normalizacion-pushforward-punto} becomes
\begin{equation}
(i_N)_!^{-1}(c)
=(-1)^N
\left\langle\operatorname{ch}_c(c),
[\mathbb C^N]_{\mathrm{BM}}\right\rangle.
\label{eq:inversa-punto-forma-cohomologica}
\end{equation}
Naturality of cohomological Thom and extension by zero identifies
the pairing in
\eqref{eq:chern-pushforward-con-signo-normal} with the pairing over
$TM$. The two sign factors are
\[
(-1)^r(-1)^N=(-1)^{N-n+N}=(-1)^n.
\]
Consequently,
\[
(i_N)_!^{-1}(j_!a)
=(-1)^n
\left\langle
\operatorname{ch}_c(a)\smile
\pi_{TM}^*\operatorname{Td}(TM\otimes_{\mathbb R}\mathbb C),
[TM]_{\mathrm{BM}}
\right\rangle.
\]
Theorem~\ref{teo:atiyah-singer-k-cerradas} proves
\eqref{eq:formula-cohomologica-indice-tm}.

If $M$ is oriented, the complex orientation of $TM$ orders each
horizontal direction first and its vertical direction next. The tangent
bundle orientation orders the $n$ horizontal directions first, followed by the
$n$ vertical directions. Passing from one ordering to the other requires
\[
(n-1)+(n-2)+\cdots+1=\frac{n(n-1)}{2}
\]
transpositions. Thus, for $u\in H_c^{2n}(TM;\mathbb Q)$,
\[
\langle u,[TM]_{\mathrm{BM}}\rangle
=(-1)^{\frac{n(n-1)}{2}}
\int_M\Phi_{TM}^{-1}(u).
\]
Multiplying this sign by $(-1)^n$ gives
$(-1)^{\frac{n(n+1)}{2}}$ and proves
\eqref{eq:formula-cohomologica-indice-base}.
\end{proof}

\section{Three geometric operators}

The following examples verify the analytic hypotheses and show how
the symbol class recovers familiar invariants.

\subsection{The Euler--de Rham operator}

First suppose that $(M,\mathbf{g})$ is a closed Riemannian manifold. Define the total odd
operator
\[
D_{\mathrm{dR}}:=d+d^*
\quad\text{on }\Omega^\bullet(M;\mathbb C).
\]
With the even and odd grading of complex forms, its restriction is
\[
D_{\mathrm{dR}}^+\colon
\Omega^{\mathrm{even}}(M;\mathbb C)
\longrightarrow
\Omega^{\mathrm{odd}}(M;\mathbb C).
\]
For $\xi\in T_x^*M\setminus\{0\}$, its two symbols are
\[
\boldsymbol{\sigma}_1(D_{\mathrm{dR}}^+)(x,\xi)=c(\xi),
\qquad
\boldsymbol{\sigma}_1^\Psi(D_{\mathrm{dR}}^+)(x,\xi)=i\,c(\xi),
\qquad
c(\xi):=\varepsilon(\xi)-\iota(\xi^\sharp),
\qquad
c(\xi)^2=-|\xi|_{\mathbf{g}}^2I.
\]
Consequently, $D_{\mathrm{dR}}^+$ is elliptic.

\begin{proposition}[Euler--de Rham index]
\label{prop:indice-euler-de-rham}
Let $(M,\mathbf{g})$ be a closed Riemannian manifold of
dimension $n$, and let $D_{\mathrm{dR}}^+$ be the Euler--de Rham operator
defined above.
Then
\begin{equation}
\operatorname{ind}_{\mathrm a}(D_{\mathrm{dR}}^+)
=\sum_{q=0}^n(-1)^q\dim_{\mathbb C}H_{\mathrm{dR}}^q(M;\mathbb C)
=\chi(M).
\label{eq:indice-euler-de-rham}
\end{equation}
If $M$ is oriented, the cohomological formula reduces to
\begin{equation}
\operatorname{ind}_{\mathrm a}(D_{\mathrm{dR}}^+)
=\int_M e(TM),
\label{eq:gauss-bonnet-indice}
\end{equation}
where $e(TM)$ is the Euler class represented in de Rham cohomology.
\end{proposition}

The identity of classes used in the proof can be found in
\cite[equation~(17.148) and Theorem~17.68, pp.~611--612]{BleeckerBoossIndex};
its integral is identified with $\chi(M)$ in
\cite[Theorem~13.6(d), pp.~320--321]{BleeckerBoossIndex}.

\begin{proof}
The square of $D_{\mathrm{dR}}=d+d^*$ is the Hodge Laplacian
$\Delta=dd^*+d^*d$. The elliptic theory already developed gives
\[
L^2\Omega^q(M;\mathbb C)
=\ker\Delta_q\oplus\operatorname{Ran}\Delta_q
\]
and a Green operator $G_q$, defined as $\Delta_q^{-1}$ on
$(\ker\Delta_q)^\perp$ and zero on $\ker\Delta_q$, which gains two
derivatives. If $\boldsymbol{\omega}\in\Omega^q(M;\mathbb C)$ and
$d\boldsymbol{\omega}=0$, write
\[
\boldsymbol{\omega}=h+\Delta G_q\boldsymbol{\omega},
\qquad h\in\ker\Delta_q.
\]
The identities $d\Delta=\Delta d$ and uniqueness of the inverse on
$(\ker\Delta)^{\perp}$ imply that $dG_q\boldsymbol{\omega}$ is harmonic. Since it is also
exact,
\[
\|dG_q\boldsymbol{\omega}\|_{L^2\Omega^{q+1}(M;\mathbb C)}^2
=\langle G_q\boldsymbol{\omega},d^*dG_q\boldsymbol{\omega}\rangle_{
L^2\Omega^q(M;\mathbb C)}=0.
\]
Therefore, $dG_q\boldsymbol{\omega}=0$ and
\[
\boldsymbol{\omega}=h+d(d^*G_q\boldsymbol{\omega}).
\]
Every de Rham cohomology class thus has a harmonic representative.
If a harmonic form is exact, the same calculation shows that it is zero;
the representative is unique. Consequently,
\[
\ker\Delta_q\cong H_{\mathrm{dR}}^q(M;\mathbb C).
\]
The kernels of $D_{\mathrm{dR}}^+$ and its adjoint are the even and odd harmonic
forms, respectively, proving the first equality in
\eqref{eq:indice-euler-de-rham}.
Theorem~\ref{teo:de-rham-euler-apendice} identifies the alternating sum of
these dimensions with $\chi(M)$.

Finally, the preceding Clifford symbol defines the \emph{Euler symbol
class} in $K_c^0(T^*M)$; it is not identified here with a Thom
class of a real bundle. Denote by
$\sharp_g\colon T^*M\to TM$ the metric isomorphism and by
$(\sharp_g)_*$ the induced isomorphism in compactly supported $K$--theory.
Its characteristic identity is
\[
(-1)^{\frac{n(n+1)}{2}}
\Phi_{TM}^{-1}\!\left(
\operatorname{ch}_c\bigl(
(\sharp_g)_*[\boldsymbol{\sigma}(D_{\mathrm{dR}}^+)]\bigr)
\right)
\operatorname{Td}(TM\otimes_{\mathbb R}\mathbb C)
=e(TM).
\]
This is an identity for the full symbol class, including the
Thom isomorphism and its orientations; it does not follow by identifying the
exterior difference with a single difference of half-spinor bundles. Indeed,
for a pair of roots $\pm x$ we have
\[
\operatorname{ch}(\Lambda^{\mathrm{even}}-\Lambda^{\mathrm{odd}})
=2-e^x-e^{-x}
=-\bigl(e^{x/2}-e^{-x/2}\bigr)^2,
\]
and the remaining factors come precisely from Thom.
Substituting this identity into
\eqref{eq:formula-cohomologica-indice-base} yields
\eqref{eq:gauss-bonnet-indice}.
\end{proof}

\subsection{The Dolbeault operator}

Now let $M$ be a closed complex manifold of complex dimension $m$ and
let $\mathbf{E}\longrightarrow M$ be a holomorphic vector bundle equipped with a
Hermitian metric. Also fix a Hermitian metric on $M$. The
complex structure splits the complexified bundles
as
\[
TM\otimes_{\mathbb R}\mathbb C=T^{1,0}M\oplus T^{0,1}M,
\qquad
T^*M\otimes_{\mathbb R}\mathbb C=T^{*1,0}M\oplus T^{*0,1}M,
\]
and we write
\[
\Omega^{p,q}(M,\mathbf{E})
=\Gamma\!\left(M,
\Lambda^pT^{*1,0}M\otimes
\Lambda^qT^{*0,1}M\otimes \mathbf{E}\right).
\]
In a holomorphic trivialization of $\mathbf{E}$ and holomorphic coordinates
$z^1,\ldots,z^m$, the twisted Dolbeault operator is
\begin{equation}
\bar\partial_{\mathbf{E}}
=\sum_{j=1}^m d\bar z^j\wedge
\frac{\partial}{\partial\bar z^j}
\colon\Omega^{p,q}(M,\mathbf{E})\longrightarrow\Omega^{p,q+1}(M,\mathbf{E}).
\label{eq:def-dolbeault-torcido-local}
\end{equation}
The transition functions of $\mathbf{E}$ are holomorphic, so
\eqref{eq:def-dolbeault-torcido-local} is independent of the trivialization.
Commutation of the derivatives $\partial/\partial\bar z^j$ and
anticommutation of exterior products prove
$\bar\partial_{\mathbf{E}}^2=0$.

The Hermitian metrics on $M$ and $\mathbf{E}$ define the formal adjoint
$\bar\partial_{\mathbf{E}}^*$. The total odd operator
\[
D_{\bar\partial,\mathbf{E}}
:=\sqrt{2}\,(\bar\partial_{\mathbf{E}}+\bar\partial_{\mathbf{E}}^*)
\quad\text{on }\Omega^{0,\bullet}(M,\mathbf{E})
\]
has as its component from even to odd degree the restriction
\[
D_{\bar\partial,\mathbf{E}}^+
:=D_{\bar\partial,\mathbf{E}}\big|_{\Omega^{0,\mathrm{even}}}
\colon
\Omega^{0,\mathrm{even}}(M,\mathbf{E})
\longrightarrow
\Omega^{0,\mathrm{odd}}(M,\mathbf{E})
\]
This operator is elliptic. Indeed, if $\xi$ is a real covariable,
the preceding local
formula and the adjoint symbol formula give, respectively,
\begin{equation}
\boldsymbol{\sigma}_1(D_{\bar\partial,\mathbf{E}}^+)(x,\xi)
=\sqrt{2}\left(
\varepsilon(\xi^{0,1})
-\varepsilon(\xi^{0,1})^*
\right)\otimes I_{\mathbf{E}_x},
\qquad
\boldsymbol{\sigma}_1^\Psi(D_{\bar\partial,\mathbf{E}}^+)(x,\xi)
=i\sqrt{2}\left(
\varepsilon(\xi^{0,1})
-\varepsilon(\xi^{0,1})^*
\right)\otimes I_{\mathbf{E}_x}.
\label{eq:simbolo-dolbeault-local}
\end{equation}
Since $2|\xi^{0,1}|_{\mathbf{g}}^2=|\xi|_{\mathbf{g}}^2$ for a real covariable and
$\varepsilon\varepsilon^*+\varepsilon^*\varepsilon
=|\xi^{0,1}|_{\mathbf{g}}^2I$, we obtain
\[
\boldsymbol{\sigma}_1^\Psi(D_{\bar\partial,\mathbf{E}}^+)(x,\xi)^*
\boldsymbol{\sigma}_1^\Psi(D_{\bar\partial,\mathbf{E}}^+)(x,\xi)=|\xi|_{\mathbf{g}}^2I.
\]
The factor $\sqrt{2}$ does not affect the symbol class. Write
\[
H^{0,q}_{\bar\partial}(M,\mathbf{E})
:=
\frac{
\ker(\bar\partial_{\mathbf{E}}\colon\Omega^{0,q}(M,\mathbf{E})\longrightarrow
\Omega^{0,q+1}(M,\mathbf{E}))
}{
\operatorname{Ran}(\bar\partial_{\mathbf{E}}\colon\Omega^{0,q-1}(M,\mathbf{E})
\longrightarrow\Omega^{0,q}(M,\mathbf{E}))
}.
\]

\begin{proposition}[Analytic Hirzebruch--Riemann--Roch]
\label{prop:indice-dolbeault}
Let $M$ be a closed complex manifold of complex
dimension $m$ equipped with a Hermitian metric; let
$\mathbf{E}\to M$ be a holomorphic vector bundle equipped with a Hermitian bundle metric, and let
$D_{\bar\partial,\mathbf{E}}^+$ be the twisted Dolbeault operator defined
above. Then
\begin{equation}
\operatorname{ind}_{\mathrm a}(D_{\bar\partial,\mathbf{E}}^+)
=\sum_{q=0}^{m}(-1)^q
\dim_{\mathbb C}H^{0,q}_{\bar\partial}(M,\mathbf{E})
=\int_M
\left[
\operatorname{ch}(\mathbf{E})\smile\operatorname{Td}(T^{1,0}M)
\right]_{2m}.
\label{eq:indice-dolbeault-rr}
\end{equation}
\end{proposition}

\begin{proof}
The square of $D_{\bar\partial,\mathbf{E}}$ is twice the Dolbeault
Laplacian,
\[
\Delta_{\bar\partial,\mathbf{E}}
=\bar\partial_{\mathbf{E}}\bar\partial_{\mathbf{E}}^*
+\bar\partial_{\mathbf{E}}^*\bar\partial_{\mathbf{E}}.
\]
It is elliptic and self-adjoint, so its range is closed and it has a
Green operator $G_{\bar\partial,\mathbf{E}}$ equal to the inverse on the orthogonal
complement of its kernel and zero on the kernel. If
$\boldsymbol{\alpha}\in\Omega^{0,q}(M,\mathbf{E})$ satisfies
$\bar\partial_{\mathbf{E}}\boldsymbol{\alpha}=0$, write
\[
\boldsymbol{\alpha}=h+\Delta_{\bar\partial,\mathbf{E}}G_{\bar\partial,\mathbf{E}}\boldsymbol{\alpha},
\qquad
h\in\ker\Delta_{\bar\partial,\mathbf{E}}.
\]
The identity
$\bar\partial_{\mathbf{E}}\Delta_{\bar\partial,\mathbf{E}}
=\Delta_{\bar\partial,\mathbf{E}}\bar\partial_{\mathbf{E}}$ shows that
$\bar\partial_EG_{\bar\partial,\mathbf{E}}\boldsymbol{\alpha}$ is harmonic. Since it is exact, its
squared norm is
\[
\left\|\bar\partial_EG_{\bar\partial,\mathbf{E}}\boldsymbol{\alpha}\right\|_{
L^2\Omega^{0,q+1}(M,\mathbf{E})}^2
=
\left\langle
G_{\bar\partial,\mathbf{E}}\boldsymbol{\alpha},
\bar\partial_{\mathbf{E}}^*\bar\partial_EG_{\bar\partial,\mathbf{E}}\boldsymbol{\alpha}
\right\rangle_{L^2\Omega^{0,q}(M,\mathbf{E})}
=0.
\]
Therefore,
\[
\boldsymbol{\alpha}
=h+\bar\partial_{\mathbf{E}}
\left(\bar\partial_{\mathbf{E}}^*G_{\bar\partial,\mathbf{E}}\boldsymbol{\alpha}\right).
\]
The same identity proves that an exact harmonic form is zero. We thus obtain,
without losing or adding any classes,
\[
\ker\Delta_{\bar\partial,\mathbf{E}}^{0,q}
\cong H^{0,q}_{\bar\partial}(M,\mathbf{E}).
\]
This gives the first equality in
\eqref{eq:indice-dolbeault-rr}.

The Clifford module
\[
\mathbf{S}_{\mathbb C}=\Lambda^{0,*}T^*M
=\Lambda^{0,\mathrm{even}}T^*M
 \oplus\Lambda^{0,\mathrm{odd}}T^*M
\]
defines the canonical $\operatorname{Spin}^c$ structure of a complex
manifold and its corresponding orientation. Formula
\eqref{eq:simbolo-dolbeault-local} shows that the symbol class is the
$\operatorname{Spin}^c$ Thom class of this module multiplied by
$\pi^*[\mathbf{E}]$; it does not depend on choosing the summand $T^{0,1}M$ in isolation.
If $x_1,\ldots,x_m$ are the formal Chern roots of $T^{1,0}M$, the
identity
\eqref{eq:defecto-chern-thom} cancels, root by root, the factor coming
from $TM\otimes_{\mathbb R}\mathbb C$ and leaves
\[
\prod_{j=1}^{m}\frac{x_j}{1-e^{-x_j}}
=\operatorname{Td}(T^{1,0}M).
\]
The complex manifold carries its canonical orientation.
Formula~\eqref{eq:formula-cohomologica-indice-base} gives the second equality
in \eqref{eq:indice-dolbeault-rr}.
\end{proof}

\subsection{The twisted Dirac operator}

Let $(M^{2m},\mathbf{g})$ be a closed oriented Riemannian spin manifold, and let
$\mathbf{E}\to M$ be a complex vector bundle equipped with a Hermitian bundle metric
$\mathbf{h}_{\mathbf{E}}$ and a compatible connection $\nabla^{\mathbf{E}}$. Clifford
geometry, the spin structure, the spinor bundle, and its connection were
constructed in Section~\ref{sec:clifford-spin-dirac}. The twisted differential
operator was defined later, once the general theory of
operators on bundles was available, in
Definition~\ref{def:dirac-torcido-geometrico}. In particular,
\[
\mathbf{S}:=\boldsymbol\Sigma M=\mathbf{S}^+\oplus \mathbf{S}^-,
\qquad
D_{\mathbf{E}}=
\begin{pmatrix}
0&D_{\mathbf{E}}^-\\
D_{\mathbf{E}}^+&0
\end{pmatrix},
\qquad
(D_{\mathbf{E}}^+)^*=D_{\mathbf{E}}^-.
\]
\[
D_{\mathbf{E}}^2=
\begin{pmatrix}
D_{\mathbf{E}}^-D_{\mathbf{E}}^+&0\\
0&D_{\mathbf{E}}^+D_{\mathbf{E}}^-
\end{pmatrix}.
\]
Both matrices are written with respect to the decomposition
$(\mathbf{S}^+\otimes\mathbf{E})\oplus
(\mathbf{S}^-\otimes\mathbf{E})$ of the bundle $\mathbf{S}\otimes\mathbf{E}$.
The corresponding Weitzenböck formula is
\eqref{eq:weitzenbock-dirac-torcido}. To calculate the index, we also use
the two symbol conventions introduced at different stages:
by Proposition~\ref{prop:simbolo-elipticidad-dirac} and
Example~\ref{ej:dirac-como-pseudodiferencial}, the differential and
Fourier symbols are, respectively,
\begin{equation}
\boldsymbol{\sigma}_1(D_{\mathbf{E}}^+)(x,\xi)
=c(\xi)\otimes I_{\mathbf{E}_x},
\qquad
\boldsymbol{\sigma}_1^\Psi(D_{\mathbf{E}}^+)(x,\xi)
=i\,c(\xi)\otimes I_{\mathbf{E}_x}.
\label{eq:dos-simbolos-dirac-torcido}
\end{equation}
The Clifford relation proves their invertibility for $\xi\neq0$.

The McKean--Singer identity expresses the Dirac index as the
supertrace of the heat semigroup associated with the square of the operator.

\begin{corollary}[McKean--Singer for the Dirac operator]
\label{cor:mckean-singer-dirac-torcido}
Let $(M^{2m},\mathbf{g})$ be a closed oriented Riemannian spin manifold, let $\mathbf{E}\to M$ be a complex vector bundle equipped with a Hermitian bundle metric
$\mathbf{h}_{\mathbf{E}}$ and a compatible connection $\nabla^{\mathbf{E}}$, and let $D_{\mathbf{E}}^+$ be the twisted Dirac operator. For every
$t>0$,
\begin{align}
\operatorname{ind}_{\mathrm a}(D_{\mathbf{E}}^+)
&=\operatorname{Tr}\bigl(e^{-tD_{\mathbf{E}}^-D_{\mathbf{E}}^+}\bigr)
-\operatorname{Tr}\bigl(e^{-tD_{\mathbf{E}}^+D_{\mathbf{E}}^-}\bigr)
\nonumber\\
&=\operatorname{Str}\bigl(e^{-tD_{\mathbf{E}}^2}\bigr)
=\int_M
\operatorname{str}_{\mathbf{S}_x\otimes \mathbf{E}_x}
\mathbf{K}_{D_{\mathbf{E}}^2}(t,x,x)\,d\lambda_{\mathbf{g}}(x),
\label{eq:mckean-singer-dirac-torcido}
\end{align}
where $\mathbf{K}_{D_{\mathbf{E}}^2}$ is the heat kernel of $D_{\mathbf{E}}^2$ and
$\operatorname{str}(A)=\operatorname{tr}(A|_{\mathbf{S}^+\otimes \mathbf{E}})
-\operatorname{tr}(A|_{\mathbf{S}^-\otimes \mathbf{E}})$.
\end{corollary}

\begin{proof}
Apply Proposition~\ref{pseudo:prop-mckean-singer} to $P=D_{\mathbf{E}}^+$.
Formula~\eqref{eq:weitzenbock-dirac-torcido} shows that $D_{\mathbf{E}}^2$ is a
Laplace-type operator with a zeroth-order potential, but its heat
kernel is not the kernel of the connection Laplacian constructed in
Chapter~\ref{cap:sobolev-fraccionario-haces-nucleo-calor}. Regularity,
the trace-class property, and the diagonal formula come from
Proposition~\ref{pseudo:prop-calor-regularizacion}. No asymptotic expansion
as $t\to 0^{+}$ is used.
\end{proof}

If the formal Chern roots of
$TM\otimes_{\mathbb R}\mathbb C$ are ordered as
$x_1,-x_1,\ldots,x_m,-x_m$,
define
\begin{equation}
\widehat A(TM)
=\prod_{j=1}^{m}
\frac{\frac{x_j}{2}}{\sinh\!\left(\frac{x_j}{2}\right)}.
\label{eq:def-clase-a-circunflejo}
\end{equation}

\begin{proposition}[Index of the Dirac operator]
\label{prop:indice-dirac-spin}
Let $(M^{2m},\mathbf{g})$ be a closed oriented Riemannian spin manifold, let $\mathbf{E}\to M$ be a complex vector bundle equipped with a Hermitian bundle metric
$\mathbf{h}_{\mathbf{E}}$ and a compatible connection $\nabla^{\mathbf{E}}$, and let $D_{\mathbf{E}}^+$ be the twisted Dirac operator. Then
\begin{equation}
\operatorname{ind}_{\mathrm a}(D_{\mathbf{E}}^+)
=\int_M
\left[
\widehat A(TM)\smile\operatorname{ch}(\mathbf{E})
\right]_{2m}.
\label{eq:indice-dirac-torcido}
\end{equation}
In particular, for $\mathbf{E}=M\times\mathbb C$ we obtain the
$\widehat A(M)$ genus.
\end{proposition}

The analytic construction and spinor Thom orientation are described
in \cite[Definitions~17.22--17.23, pp.~531--532]{BleeckerBoossIndex} and
\cite[Appendix~C, Theorem~C.12; Chapter~III, \S12,
Proposition~12.5]{LawsonMichelsohn1989}.

\begin{proof}
The second formula in \eqref{eq:dos-simbolos-dirac-torcido} shows that the symbol
class is the complex
\[
\pi^*(\mathbf{S}^+\otimes \mathbf{E})
\xrightarrow{\,i c(\xi)\otimes I_{\mathbf{E}}\,}
\pi^*(\mathbf{S}^-\otimes \mathbf{E}),
\qquad \xi\in T^*M,
\]
exact away from the zero section. The complex without $\mathbf{E}$ is the spinor Thom
class, and the twisting factor is $\pi^*[\mathbf{E}]$; multiplying the morphism by $i$
does not change its class. The Chern defect of the spinor Thom class,
stated in
Theorem~\ref{teo:defecto-chern-thom-spinorial-apendice}, is
\[
\Phi_{TM}^{-1}\operatorname{ch}_c(s(TM))
=(-1)^m\widehat A(TM)^{-1}.
\]
Moreover, the splitting principle applied to the pairs of roots
$x_j,-x_j$ gives
\[
\operatorname{Td}(TM\otimes_{\mathbb R}\mathbb C)
=\widehat A(TM)^2.
\]
The factor $\pi^*[\mathbf{E}]$ contributes $\operatorname{ch}(\mathbf{E})$. Since the overall sign
in \eqref{eq:formula-cohomologica-indice-base} in dimension $2m$ is
$(-1)^m$, the two signs cancel, leaving
$\widehat A(TM)\operatorname{ch}(\mathbf{E})$. This proves
\eqref{eq:indice-dirac-torcido}.
\end{proof}

\begin{proposition}[Index of the $\operatorname{Spin}^c$ Dirac operator]
\label{prop:indice-dirac-spinc}
Let $M^{2m}$ be a closed oriented smooth manifold equipped with a
$\operatorname{Spin}^c$ structure with determinant line $\mathbf{L}$. Let $A$ be a unitary connection
on $\mathbf{L}$ and let $\mathbf{E}$ be a complex vector bundle equipped with a Hermitian bundle metric
$\mathbf{h}_{\mathbf{E}}$ and a compatible connection $\nabla^{\mathbf{E}}$. The
operator $D_{A,\mathbf{E}}$ is obtained by twisting the operator of
Proposition~\ref{def:operadores-dirac-spinc} by $\nabla^{\mathbf{E}}$, using
the same product connection construction as in
Definition~\ref{def:dirac-torcido-geometrico}. The geometric data of the
$\operatorname{Spin}^c$ structure are those of
Section~\ref{sec:clifford-spin-dirac}. Its chiral restriction satisfies
\[
D_{A,\mathbf{E}}^+\colon
\Gamma(\Sigma^{c,+}M\otimes \mathbf{E})
\longrightarrow
\Gamma(\Sigma^{c,-}M\otimes \mathbf{E}).
\]
Then
\begin{equation}
\operatorname{ind}_{\mathrm a}(D_{A,\mathbf{E}}^+)
=\int_M
\left[
\widehat A(TM)\smile e^{c_1(\mathbf{L})/2}
\smile\operatorname{ch}(\mathbf{E})
\right]_{2m}.
\label{eq:indice-dirac-spinc}
\end{equation}
If the structure comes from an almost complex structure, then
$\mathbf{L}=\det_{\mathbb C}(T^{1,0}M)$ and
\[
\widehat A(TM)e^{c_1(\mathbf{L})/2}
=\operatorname{Td}(T^{1,0}M).
\]
If the almost complex structure is Kähler, the connection on
$\mathbf{L}=K_M^{-1}$ is the one induced by Levi--Civita, and $\mathbf{E}$ is holomorphic with
its Chern connection, the canonical operator is exactly
$\sqrt2(\bar\partial_{\mathbf{E}}+\bar\partial_{\mathbf{E}}^*)$. In the general almost complex
case, the $\operatorname{Spin}^c$ Dirac operator and the associated Dolbeault-type
operator have the same symbol and index, although they may differ
by a zeroth-order term.
\end{proposition}

The canonical topological identifications are collected in
\cite[Appendix~D, Theorem~D.2 and Examples~D.5--D.6;
Appendix~C, Theorem~C.12; Appendix~D,
Theorem~D.15]{LawsonMichelsohn1989}.

\begin{proof}
The $\operatorname{Spin}^c$ Thom class is obtained from the spinor class
formally twisted by the square root of $\mathbf{L}$; in the Chern character this
introduces the well-defined factor $e^{c_1(\mathbf{L})/2}$. The same Thom
calculation and \eqref{eq:formula-cohomologica-indice-base} used in the preceding
proposition give \eqref{eq:indice-dirac-spinc}. For the almost complex structure,
the determinant line and characteristic identity follow, under the
splitting principle, from
\[
\frac{x/2}{\sinh(x/2)}e^{x/2}
=\frac{x}{1-e^{-x}}.
\]
The Kähler identification follows from the Chern--Levi--Civita connection
on $\Lambda^{0,\bullet}T^*M$.
\end{proof}

These cases illustrate different roles of the same theorem. For de Rham,
the index is the Euler characteristic; for Dolbeault, the holomorphic
Euler characteristic; in the spinor case, a characteristic number built from
$\widehat A(TM)$. The equality is independent of the auxiliary metrics and
connections: smooth variation preserves the symbol class and gives a
continuous family of Fredholm operators.

\part{Geometric evolution equations}
\chapter*{Introduction to Part V}
\addcontentsline{toc}{chapter}{Introduction to Part V}
\markboth{Part V. Geometric evolution equations}{Introduction to Part V}

Until now, the Riemannian metric has been a fixed part of the geometry. In this part we allow it to change with time. The idea is simple to state: instead of seeking a metric with special properties directly, we start with an initial metric and let it evolve according to an equation determined by its curvature. The geometric problem thereby becomes a nonlinear partial differential equation.

Before studying Ricci flow, we need a parabolic theory that establishes existence, uniqueness, and regularity for systems on manifolds. The flow also presents a difficulty specific to geometry: the equation is invariant under diffeomorphisms, and this symmetry appears analytically as a degeneracy of the symbol. DeTurck's procedure removes this degeneracy through a choice of gauge and reduces the problem, for a short time, to a strictly parabolic system.

Once the solution has been constructed, the evolution equations for curvature and the maximum principles yield a priori estimates. When curvature grows and a singularity appears, rescaling allows the geometry to be observed at its natural scale. Compactness extracts limits of these magnified regions, while Perelman's monotone quantities and noncollapsing restrict the possible models. In dimension three, this analysis leads to surgery and shows particularly clearly how an evolution equation can reveal the topological structure of a manifold.

\begin{semblanzaHistorica}{Hamilton, DeTurck, and Perelman}
Hamilton introduced Ricci flow in 1982 with the idea of deforming a metric using its own curvature. It soon became clear that diffeomorphism invariance required careful treatment; DeTurck's procedure supplied the adjustment needed to apply parabolic theory. Two decades later, Perelman introduced new monotone quantities, the noncollapsing theorem, and a decisive analysis of singularities and surgery. This history illustrates a characteristic feature of geometric analysis: an intuition about the shape of space leads to an equation, the equation demands precise analytic tools, and those tools ultimately return to a geometric or topological problem.
\end{semblanzaHistorica}

\chapter{Geometric parabolic equations on closed manifolds}
\label{cap:teoria-parabolica-cerradas}

A geometric evolution equation resembles the heat equation only
after several difficulties have been resolved. The unknown is usually a
section of a bundle, the coefficients of the principal part depend on that
same section, and geometric symmetries may conceal the parabolicity of the
system. Before studying Ricci flow, it is useful to separate these problems
and develop a theory that identifies precisely which part of the
argument comes from diffusion and which depends on the nonlinear structure.

We work on a closed Riemannian manifold $(M,\mathbf{g}_0)$ of
dimension $m$. Let $\mathbf{E}\to M$ be a real vector bundle equipped with a
bundle metric $\mathbf{h}_{\mathbf{E}}$ and a compatible connection
$\nabla^{\mathbf{E}}$, abbreviated to $\nabla$ together with its induced
connections. These structures are used to measure sections and write the
equations; any other smooth choice gives equivalent norms.
Denote by $H^r(M,\mathbf{E})$ the Sobolev space of integer order $r$
constructed in Chapter~\ref{cap:sobolev-haces}. The operator
$\Delta_\nabla=\nabla^*\nabla$ defining the spectral scale is nonnegative,
whereas the covariant trace $\operatorname{tr}_{\mathbf{g}_0}(\nabla^2)$
used in the maximum principles has the opposite sign.

Throughout this chapter, fix an integer
\begin{equation}
\label{eq:orden-sobolev-parabolico}
 s>\frac{m}{2}+2.
\end{equation}
Set $\Lambda=(I+\Delta_\nabla)^{\displaystyle\frac{1}{2}}$. By
Corollary~\ref{cor:bessel-sobolev-entero-variedad-cerrada}, the norm
\begin{equation}
\label{eq:norma-espectral-parabolica}
 \|\mathbf{u}\|_{H^r(M,\mathbf{E})}:=\|\Lambda^r\mathbf{u}\|_{L^2(M,\mathbf{E})}
\end{equation}
is equivalent to the covariant norm when $r\geq0$ is an integer; for $r<0$
it is interpreted by completion. We use this choice because it makes
energy identities in the Sobolev scale exact. Estimates
obtained with it remain valid for the covariant norms
by equivalence. In particular, $H^{s-2}(M)$, $H^{s-1}(M)$, and
$H^s(M)$ are algebras and the inclusion
$H^{s-1}(M,\mathbf{E})\hookrightarrow C^0(M,\mathbf{E})$ is continuous.

The linear estimate must control both the solution and its time
derivative and retain its constants as the coefficients vary within
fixed bounds. This uniformity allows us to solve the quasilinear equation:
each approximation determines the coefficients of the next, while
a common energy controls the entire sequence. Maximum principles
will then provide pointwise information that integral norms alone
do not supply.

\begin{semblanzaHistorica}{From the heat equation to geometric flows}
Parabolic theory grew out of the study of diffusion: heat dampens
irregularities and produces smoothness at positive times. In a geometric
flow the unknown may be a metric, a connection, or a section, and
the coefficients of the equation change with it. Energy estimates,
maximum principles, and product inequalities preserve the
smoothing mechanism of heat in this nonlinear regime. The
results of this chapter prepare this transition before applying it to
Ricci flow.
\end{semblanzaHistorica}

\section{Product estimates and energy spaces}

The nonlinear estimates we need are precise consequences of
the Leibniz rule and the Gagliardo--Nirenberg inequality in
Theorem~\ref{teo:gagliardo-nirenberg-haces}. We first collect them to make
explicit where condition~\eqref{eq:orden-sobolev-parabolico} is used.

\begin{lemma}[Sobolev product with distributed orders]
\label{lem:producto-sobolev-ordenes-distribuidos-parabolico}
Let $\mathbf{V}_1,\mathbf{V}_2,\mathbf{V}_3\to M$ be vector bundles with
fixed metrics and connections, and let
$\beta\colon\mathbf{V}_1\otimes\mathbf{V}_2\to\mathbf{V}_3$ be a smooth
homomorphism. If $a,b\in\mathbb N_0$ and $a+b>\displaystyle\frac{m}{2}$, there exists a constant $C$ such that
\begin{equation}
 \|\beta(\mathbf{u},\mathbf{v})\|_{L^2(M,\mathbf{V}_3)}
 \leq C\|\mathbf{u}\|_{H^a(M,\mathbf{V}_1)}
       \|\mathbf{v}\|_{H^b(M,\mathbf{V}_2)}.
 \label{eq:producto-sobolev-ordenes-distribuidos-parabolico}
\end{equation}
The same estimate holds if one factor is a covariant derivative of
a section and the orders $a,b$ describe the regularity remaining
after taking that derivative.
\end{lemma}

\begin{proof}
By a finite partition of unity and local equivalence of Sobolev
norms, it suffices to prove the assertion for vector-valued functions supported in
a relatively compact Euclidean open subset. Since $a+b>\displaystyle\frac{m}{2}$, we can choose
$p,q\in[2,\infty]$ so that
\[
 \frac{1}{p}+\frac{1}{q}=\frac{1}{2},
 \qquad H^a(\mathbb R^m)\hookrightarrow L^p(\mathbb R^m),
 \qquad H^b(\mathbb R^m)\hookrightarrow L^q(\mathbb R^m).
\]
Let us specify the choice in the critical cases. If $a>\displaystyle\frac{m}{2}$, take
$p=\infty$ and $q=2$; if $b>\displaystyle\frac{m}{2}$, proceed symmetrically. Thus suppose
that $a,b\leq \displaystyle\frac{m}{2}$. When both are less than $\displaystyle\frac{m}{2}$, the inequality
$a+b>\displaystyle\frac{m}{2}$ allows us to choose a number $\theta$ such that
\[
 \frac{1}{2}-\frac{a}{m}<\theta<\frac{b}{m},
 \qquad \frac{1}{p}=\theta,
 \qquad \frac{1}{q}=\frac{1}{2}-\theta.
\]
The corresponding Sobolev inclusions are continuous. If $a=\displaystyle\frac{m}{2}$,
choose $0<\theta<\displaystyle\frac{b}{m}$; the inclusion
$H^{\displaystyle\frac{m}{2}}(\mathbb R^m)\hookrightarrow L^p(\mathbb R^m)$ holds for every $p<\infty$. The case
$b=\displaystyle\frac{m}{2}$ is analogous, taking $\displaystyle\frac{1}{2}-\displaystyle\frac{a}{m}<\theta<\displaystyle\frac{1}{2}$. This also
covers the critical orders using finite exponents. The
inclusion $H^{\displaystyle\frac{m}{2}}(\mathbb R^m)\hookrightarrow L^\infty(\mathbb R^m)$ does not hold in general. Hölder and the two inclusions give
\[
 \|\beta(\mathbf{u},\mathbf{v})\|_{L^2(M,\mathbf V_3)}
 \leq C\|\mathbf{u}\|_{L^p(M,\mathbf V_1)}\|\mathbf{v}\|_{L^q(M,\mathbf V_2)}
 \leq C\|\mathbf{u}\|_{H^a(M,\mathbf V_1)}\|\mathbf{v}\|_{H^b(M,\mathbf V_2)}.
\]
The coefficients of $\beta$, the connections, and the local frames are
smooth on finitely many compact sets and can therefore be absorbed into the
constant.
\end{proof}

\begin{lemma}[Moser estimates]
\label{lem:moser-parabolico-haces}
Let $(M,\mathbf{g}_0)$ be a closed Riemannian manifold of
dimension $m$.
In the following estimates, for a section of a bundle
$\mathbf V\to M$, write $\|\cdot\|_\infty$ for the norm of
$L^\infty(M,\mathbf V)$ with the fixed metrics.
Let $\mathbf{V}_1,\mathbf{V}_2,\mathbf{V}_3\longrightarrow M$ be vector bundles with fixed metrics and
connections, and let
$\beta\colon \mathbf{V}_1\otimes \mathbf{V}_2\longrightarrow \mathbf{V}_3$ be a smooth homomorphism.
For every integer $r\geq0$ there exists a constant $C$, depending on the
fixed structures and $r$, such that, for
$\mathbf{u}\in H^r(M,\mathbf{V}_1)\cap L^\infty(M,\mathbf{V}_1)$ and
$\mathbf{v}\in H^r(M,\mathbf{V}_2)\cap L^\infty(M,\mathbf{V}_2)$,
\begin{align}
 \|\beta(\mathbf{u},\mathbf{v})\|_{H^r(M,\mathbf{V}_3)}
 &\leq C\bigl(
 \|\mathbf{u}\|_\infty\|\mathbf{v}\|_{H^r(M,\mathbf{V}_2)}
 +\|\mathbf{u}\|_{H^r(M,\mathbf{V}_1)}\|\mathbf{v}\|_\infty
 \bigr),
 \label{eq:moser-producto-parabolico}\\
 \|[\nabla^r,\mathbf{a}]\mathbf{v}\|_{L^2(M,T^{(0,r)}(TM)\otimes \mathbf{V}_3)}
 &\leq C\bigl(
 \|\nabla \mathbf{a}\|_\infty
 \|\mathbf{v}\|_{H^{r-1}(M,\mathbf{V}_2)}
 +\|\mathbf{a}\|_{H^r(M,\operatorname{Hom}(\mathbf{V}_2,\mathbf{V}_3))}
 \|\mathbf{v}\|_\infty
 \bigr)
 \label{eq:moser-conmutador-parabolico}
\end{align}
In the commutator estimate, assume $r\geq1$,
$\mathbf{v}\in H^{r-1}(M,\mathbf{V}_2)\cap L^\infty(M,\mathbf{V}_2)$, and
$\mathbf{a}\in H^r(M,\operatorname{Hom}(\mathbf{V}_2,\mathbf{V}_3))$ with
$\nabla \mathbf{a}\in L^\infty(M,T^*M\otimes\operatorname{Hom}(\mathbf{V}_2,\mathbf{V}_3))$. The
action of $\mathbf{a}$ on the tensor factors is the identity on
$T^{(0,r)}(TM)$, and
\[
 [\nabla^r,\mathbf{a}]\mathbf{v}
 :=\nabla^r(\mathbf{a}\mathbf{v})
 -(I_{T^{(0,r)}(TM)}\otimes \mathbf{a})(\nabla^r\mathbf{v})
 \in L^2\bigl(M,T^{(0,r)}(TM)\otimes \mathbf{V}_3\bigr).
\]
For nonsmooth sections, the difference defining the commutator is taken
first in distributions; the estimate states that this distribution is
represented by a section of $L^2(M,T^{(0,r)}(TM)\otimes\mathbf V_3)$.
If
$\mathbf{a}\in H^s(M,\operatorname{Hom}(T^*M\otimes \mathbf{V}_1,\mathbf{V}_3))$ and
$\mathbf{u}\in H^r(M,\mathbf{V}_1)$, the same commutator applied to $\nabla \mathbf{u}$ satisfies
\begin{equation}
\label{eq:moser-conmutador-asimetrico-parabolico}
 \|[\nabla^r,\mathbf{a}](\nabla \mathbf{u})\|_{L^2(M,T^{(0,r)}(TM)\otimes \mathbf{V}_3)}
 \leq C
 \|\mathbf{a}\|_{H^s(M,\operatorname{Hom}(T^*M\otimes \mathbf{V}_1,\mathbf{V}_3))}
 \|\mathbf{u}\|_{H^r(M,\mathbf{V}_1)},
 \qquad 1\leq r\leq s.
\end{equation}

Now let $\mathcal U\subseteq \mathbf{V}_1$ be an open subset of a bundle and let
$G\colon\mathcal U\longrightarrow \mathbf{V}_2$ be a smooth bundle map. If
$K\subseteq\mathcal U$ is compact, there exists a neighborhood $\mathcal V$ of $K$
and, for each $R>0$, a constant $C_{K,R}$ with the following properties.
If $\mathbf{u},\mathbf{v}\in H^r(M,\mathbf{V}_1)$, the images of $\mathbf{u}$ and $\mathbf{v}$ lie in $\mathcal V$, and
$\|\mathbf{u}\|_{H^r(M,\mathbf{V}_1)}+\|\mathbf{v}\|_{H^r(M,\mathbf{V}_1)}\leq R$, then
\begin{align}
 \|G(\mathbf{u})\|_{H^r(M,\mathbf{V}_2)}&\leq C_{K,R},
 \label{eq:moser-composicion-parabolico}\\
 \|G(\mathbf{u})-G(\mathbf{v})\|_{H^r(M,\mathbf{V}_2)}
 &\leq C_{K,R}\|\mathbf{u}-\mathbf{v}\|_{H^r(M,\mathbf{V}_1)}.
 \label{eq:moser-composicion-lipschitz-parabolico}
\end{align}
These last two estimates hold whenever $r>\displaystyle\frac{m}{2}$.
Moreover, for every integer $r\geq0$ and every section $\mathbf{u}$ with image in
$\mathcal V$, we have the higher-order estimate
\begin{equation}
 \|G(\mathbf{u})\|_{H^r(M,\mathbf{V}_2)}
 \leq C_{K,r}\bigl(1+\|\mathbf{u}\|_{H^r(M,\mathbf{V}_1)}\bigr).
 \label{eq:moser-composicion-lineal-orden-alto}
\end{equation}
Here the constant depends on a uniform bound for the values of $\mathbf{u}$
and derivatives of $G$ up to order $r$ on a fixed neighborhood of those values;
it does not depend on the order-$r$ norm of $\mathbf{u}$.
These estimates also hold for maps depending smoothly on several arguments, with
the norm of the difference of all of them on the right-hand side.
\end{lemma}

\begin{proof}
Localize using a finite partition of unity and work in a
chart and local frame. The covariant norms and Euclidean norms of the
components are equivalent by
Lemma~\ref{lem:local-expression-higher-order}. Its triangular formula shows
which connection coefficients enter at each order; on the
compact sets of this localization, all their required derivatives are bounded. It therefore suffices to justify the estimates
for partial derivatives of vector-valued functions. Fix a regular coordinate ball $(U,\phi)$ and set $\Omega:=\phi(U)\subset\mathbb R^m$. In the following local calculations, $u,v,a$ denote the scalar components of the sections or coefficients in the fixed frames; their Euclidean norms are taken over $\Omega$.

Let $1\leq k\leq r$ and $|\gamma|=k$. The Leibniz rule gives
\[
 \partial^\gamma(uv)=\sum_{\substack{\mu\in\mathbb N_0^m\\\mu\leq\gamma}}
 \binom{\gamma}{\mu}(\partial^\mu u)
 (\partial^{\gamma-\mu}v).
\]
The terms $\mu=0$ and $\mu=\gamma$ are bounded, respectively, by
$\|u\|_{L^\infty(\Omega)}\|v\|_{H^k(\Omega)}$ and
$\|u\|_{H^k(\Omega)}\|v\|_{L^\infty(\Omega)}$. If
$0<|\mu|<k$, set $\theta=\displaystyle\frac{|\mu|}{k}$ and choose
$p=\displaystyle\frac{2}{\theta}$, $q=\displaystyle\frac{2}{1-\theta}$. The Gagliardo--Nirenberg
inequality gives
\begin{align*}
 \|\partial^\mu u\|_{L^p(\Omega)}
 &\leq C\|u\|_{L^\infty(\Omega)}^{1-\theta}\|u\|_{H^k(\Omega)}^{\theta},\\
 \|\partial^{\gamma-\mu}v\|_{L^q(\Omega)}
 &\leq C\|v\|_{L^\infty(\Omega)}^{\theta}\|v\|_{H^k(\Omega)}^{1-\theta}.
\end{align*}
Since $\displaystyle\frac{1}{p}+\displaystyle\frac{1}{q}=\displaystyle\frac{1}{2}$, Hölder and the weighted Young inequality give
\begin{align*}
 \|(\partial^\mu u)(\partial^{\gamma-\mu}v)\|_{L^2(\Omega)}
 &\leq C
 \bigl(\|u\|_{H^k(\Omega)}\|v\|_{L^\infty(\Omega)}\bigr)^\theta
 \bigl(\|u\|_{L^\infty(\Omega)}\|v\|_{H^k(\Omega)}\bigr)^{1-\theta}\\
 &\leq C\bigl(
 \|u\|_{H^k(\Omega)}\|v\|_{L^\infty(\Omega)}
 +\|u\|_{L^\infty(\Omega)}\|v\|_{H^k(\Omega)}\bigr).
\end{align*}
Summing over $|\gamma|\leq r$ proves
\eqref{eq:moser-producto-parabolico}.

For the commutator, observe that
\[
 [\partial^\gamma,a]v
 =\sum_{\substack{\mu\in\mathbb N_0^m\\0<\mu\leq\gamma}}\binom{\gamma}{\mu}
 (\partial^\mu a)(\partial^{\gamma-\mu}v).
\]
Each summand distributes at most $k-1$ derivatives between $\partial a$ and $v$.
The product estimate just proved, applied at order $k-1$, gives
\[
 \|[\partial^\gamma,a]v\|_{L^2(\Omega)}
 \leq C\bigl(
 \|\partial a\|_{L^\infty(\Omega)}\|v\|_{H^{k-1}(\Omega)}
 +\|\partial a\|_{H^{k-1}(\Omega)}\|v\|_{L^\infty(\Omega)}\bigr),
\]
and summing over $k\leq r$ yields
\eqref{eq:moser-conmutador-parabolico}.

Now consider the asymmetric commutator. For $|\gamma|=k\leq r$, its
components are sums of terms
\[
 (\partial^\mu a)\partial^{\gamma-\mu+1}u,
 \qquad 1\leq|\mu|\leq k.
\]
After taking the indicated derivatives, the first factor retains
$s-|\mu|$ derivatives in $L^2(\mathbb R^m)$ and the second retains $|\mu|-1$. The sum of
these orders is $s-1>\displaystyle\frac{m}{2}$.
Lemma~\ref{lem:producto-sobolev-ordenes-distribuidos-parabolico} implies
\[
 \|(\partial^\mu a)\partial^{\gamma-\mu+1}u\|_{L^2(\Omega)}
 \leq C\|a\|_{H^s(\Omega)}\|u\|_{H^k(\Omega)}.
\]
Summing over all multiindices yields
\eqref{eq:moser-conmutador-asimetrico-parabolico}.

For the composition estimates, choose an open subset
$\mathcal V$ whose closure is compact, contains $K$, and is contained in
$\mathcal U$. A cutoff function on the total space, equal to one on
$\mathcal V$ and supported in $\mathcal U$, allows us to extend $G$ to a
smooth bundle map $\widetilde G$ defined on all fibers. The
Faà di Bruno formula of Theorem~\ref{faa di bruno multivariable}, including
derivatives in the base variable, expresses every derivative of order $k\leq r$ of
$\widetilde G(x,u(x))$ as terms of the form
\[
 (\partial_x^\alpha D_z^j\widetilde G)(x,u(x))
 [\partial^{\mu_1}u,\ldots,\partial^{\mu_j}u],
 \qquad
 \begin{gathered}
 1\leq j\leq k,\quad \alpha,\mu_1,\ldots,\mu_j\in\mathbb N_0^m,\\
 |\mu_i|\geq1\ (1\leq i\leq j),\qquad
 |\alpha|+\sum_{i=1}^j|\mu_i|=k,
 \end{gathered}
\]
and uniformly bounded terms with no derivatives of $u$. For $r\geq1$,
choose $p_i=\displaystyle\frac{2r}{|\mu_i|}$. The Gagliardo--Nirenberg inequality gives
\[
 \|\partial^{\mu_i}u\|_{L^{p_i}(\Omega)}
 \leq C\|u\|_{L^\infty(\Omega)}^{1-\frac{|\mu_i|}{r}}
          \|u\|_{H^r(\Omega)}^{\frac{|\mu_i|}{r}}.
\]
If $\displaystyle \theta:=\displaystyle\sum_{i=1}^j\displaystyle\frac{|\mu_i|}{r}$, then $0<\theta\leq1$ and
$\displaystyle \sum_{i=1}^j\displaystyle\frac{1}{p_i}=\displaystyle\frac{\theta}{2}\leq\displaystyle\frac{1}{2}$. Hölder, followed by the inclusion
$L^{\displaystyle\frac{2}{\theta}}(\Omega)\hookrightarrow L^2(\Omega)$, valid because $\Omega$ has finite measure, bounds the
product by $C_K\|u\|_{H^r(\Omega)}^{\theta}$. Indeed, for every $h\in L^{\displaystyle\frac{2}{\theta}}(\Omega)$,
\[
 \|h\|_{L^2(\Omega)}\leq |\Omega|^{\frac{1-\theta}{2}}\|h\|_{L^{\frac{2}{\theta}}(\Omega)}.
\]
The inequality
$x^\theta\leq1+x$ for $x\geq0$ proves
\eqref{eq:moser-composicion-lineal-orden-alto}. The case $r=0$ follows by
integrating the uniform bound for $G$ on the image of $u$.
The finite sum over charts preserves these estimates. In particular,
a bound $\|u\|_{H^r(M,\mathbf V_1)}\leq R$ gives
\eqref{eq:moser-composicion-parabolico}.
Finally,
\[
 \widetilde G(u)-\widetilde G(v)
 =\int_0^1D\widetilde G\bigl(v+\theta(u-v)\bigr)[u-v] \,d\theta.
\]
The global extension ensures that the segment always lies within the domain.
Applying the preceding estimate to the integrand, with constants uniform in
$\theta$, yields
\eqref{eq:moser-composicion-lipschitz-parabolico}. The same argument handles
several arguments and completes the proof.
\end{proof}

For $T>0$ and $r\geq0$ define
\begin{align*}
 \mathcal X_T^r(\mathbf{E}):={}&
 \bigl\{\mathbf{u}\in L^2(0,T;H^{r+1}(M,\mathbf{E}))\mathrel{\big|}
 \partial_t\mathbf{u}\in L^2(0,T;H^{r-1}(M,\mathbf{E}))\bigr\},\\
 \mathcal Y_T^r(\mathbf{E}):={}&L^2(0,T;H^{r-1}(M,\mathbf{E})),
\end{align*}
where the time derivative is understood in the sense of distributions with
values in $H^{r-1}(M,\mathbf{E})$. Initially equip
$\mathcal X_T^r(\mathbf{E})$ with the graph norm
\begin{equation}
 \|\mathbf{u}\|_{\mathcal X_T^r(\mathbf{E}),0}^2
 :=\int_0^T\|\mathbf{u}(t)\|_{H^{r+1}(M,\mathbf{E})}^2\,dt
 +\int_0^T\|\partial_t\mathbf{u}(t)\|_{H^{r-1}(M,\mathbf{E})}^2\,dt.
 \label{eq:norma-grafica-espacio-energia-parabolico}
\end{equation}
This norm is complete. Indeed, if $(\mathbf{u}_j)$ is Cauchy,
there exist $\mathbf{u}\in L^2(0,T;H^{r+1}(M,\mathbf E))$ and
$\mathbf{v}\in L^2(0,T;H^{r-1}(M,\mathbf E))$ such that
$\mathbf{u}_j\to\mathbf{u}$ and $\partial_t\mathbf{u}_j\to\mathbf{v}$ in the
corresponding spaces. Passing to the limit in the identity defining
the distributional derivative gives $\partial_t\mathbf{u}=\mathbf{v}$.
The next lemma allows the supremum in time to be incorporated into the norm.

We use the pairing of the Hilbert triple
$H^{r+1}(M,\mathbf{E})\subseteq H^r(M,\mathbf{E})\subseteq H^{r-1}(M,\mathbf{E})$ given, on
smooth sections, by
\begin{equation}
\label{eq:emparejamiento-terna-parabolica}
 \langle \mathbf{f},\mathbf{v}\rangle_{H^{r-1}(M,\mathbf E),H^{r+1}(M,\mathbf E)}
 :=(\Lambda^{r-1}\mathbf{f},\Lambda^{r+1}\mathbf{v})_{L^2(M,\mathbf{E})}.
\end{equation}
Thus $H^{r-1}(M,\mathbf E)$ plays the role of the dual of $H^{r+1}(M,\mathbf E)$ when $H^r(M,\mathbf E)$ is the
pivot space.

\begin{lemma}[Time trace]
\label{lem:traza-temporal-parabolica}
Let $(M,\mathbf{g}_0)$ be a closed Riemannian manifold and
let $\mathbf{E}\to M$ be a vector bundle with a bundle metric and compatible
connection.
Every element of $\mathcal X_T^r(\mathbf{E})$ has a unique representative in
$C([0,T];H^r(M,\mathbf{E}))$. For this representative,
\begin{equation}
\label{eq:identidad-energia-traza-temporal}
 \|\mathbf{u}(t)\|_{H^r(M,\mathbf{E})}^2-\|\mathbf{u}(\tau)\|_{H^r(M,\mathbf{E})}^2
 =2\int_\tau^t
 \langle\partial_\rho \mathbf{u}(\rho),\mathbf{u}(\rho)\rangle_{H^{r-1}(M,\mathbf E),H^{r+1}(M,\mathbf E)}\,d\rho
\end{equation}
for $0\leq\tau\leq t\leq T$. In particular, evaluation
$\mathbf{u}\mapsto \mathbf{u}(0)$ is continuous from $\mathcal X_T^r(\mathbf{E})$ to $H^r(M,\mathbf{E})$.
\end{lemma}

\begin{proof}
Let $\Delta_\nabla=\nabla^*\nabla$ be the connection Laplacian and set
$\Lambda=(I+\Delta_\nabla)^{\displaystyle\frac{1}{2}}$. By
Theorem~\ref{teo:espectro-laplaciano-cerrada-eigenbase}, when
\(L^2(M,\mathbf E)\) is infinite dimensional there is an orthonormal
basis \((\mathbf e_j)\) of smooth eigensections. In dimension zero,
\(M\) is a finite set and \(\Delta_\nabla=0\), so a finite basis
is used; if the rank is zero, the space is \(\{0\}\) and the basis
is empty. Let \(P_N\) be the projection onto the first \(N\) vectors.
In the case of finite dimension \(d_H\), set \(P_N=I\) for
\(N\geq d_H\).

The projections $P_N$ commute with
all powers of $\Lambda$, are contractions on each
$H^\sigma(M,\mathbf{E})$, and converge strongly to the identity on these
spaces.

Let $\mathbf{u}\in\mathcal X_T^r(\mathbf{E})$ and set
$\mathbf{u}_N=P_N\mathbf{u}$. Since
$\partial_t\mathbf{u}_N=P_N\partial_t\mathbf{u}$ in the distributional
sense, each scalar coefficient of $\mathbf{u}_N$ belongs to
$H^1(0,T)$. Thus
$\mathbf{u}_N\in C([0,T];H^r(M,\mathbf{E}))$, and the differentiation rule on the
finite-dimensional subspace gives
\begin{equation}
 \|\mathbf{u}_N(t)\|_{H^r(M,\mathbf E)}^2-\|\mathbf{u}_N(\tau)\|_{H^r(M,\mathbf E)}^2
 =2\int_\tau^t
 \langle\partial_\rho\mathbf{u}_N(\rho),\mathbf{u}_N(\rho)
 \rangle_{H^{r-1}(M,\mathbf E),H^{r+1}(M,\mathbf E)}\,d\rho.
 \label{eq:identidad-energia-galerkin-traza}
\end{equation}

Let us prove that $(\mathbf{u}_N)$ converges in
$C([0,T];H^r(M,\mathbf{E}))$. If
$\mathbf{w}=\mathbf{u}_N-\mathbf{u}_M$, there exists
$\tau\in[0,T]$ such that
\[
 \|\mathbf{w}(\tau)\|_{H^r(M,\mathbf E)}^2
 \leq\frac{1}{T}\int_0^T\|\mathbf{w}(\rho)\|_{H^r(M,\mathbf E)}^2\,d\rho
 \leq\frac{1}{T}\|\mathbf{w}\|_{L^2(0,T;H^{r+1}(M,\mathbf E))}^2.
\]
Applying \eqref{eq:identidad-energia-galerkin-traza} to $\mathbf{w}$ and
using Cauchy--Schwarz in the pairing of
\eqref{eq:emparejamiento-terna-parabolica} yields, for every
$t\in[0,T]$,
\begin{align*}
 \|\mathbf{w}(t)\|_{H^r(M,\mathbf E)}^2
 &\leq \frac{1}{T}\|\mathbf{w}\|_{L^2(0,T;H^{r+1}(M,\mathbf E))}^2
 +2\|\partial_t\mathbf{w}\|_{L^2(0,T;H^{r-1}(M,\mathbf E))}
    \|\mathbf{w}\|_{L^2(0,T;H^{r+1}(M,\mathbf E))}\\
 &\leq
 \left(1+\frac{1}{T}\right)
 \|\mathbf{w}\|_{L^2(0,T;H^{r+1}(M,\mathbf E))}^2
 +\|\partial_t\mathbf{w}\|_{L^2(0,T;H^{r-1}(M,\mathbf E))}^2.
\end{align*}
Since $P_N\to I$ strongly in both Bochner spaces occurring here,
the right-hand side tends to zero as $M,N\to\infty$. Consequently, there exists
a section
$\widetilde{\mathbf{u}}\in C([0,T];H^r(M,\mathbf{E}))$ such that
$\mathbf{u}_N\to\widetilde{\mathbf{u}}$ uniformly on $H^r(M,\mathbf E)$.
On the other hand, $\mathbf{u}_N\to\mathbf{u}$ in
$L^2(0,T;H^{r+1}(M,\mathbf E))$; hence
$\widetilde{\mathbf{u}}=\mathbf{u}$ for almost every time. Thus
$\widetilde{\mathbf{u}}$ is a representative of $\mathbf{u}$. Two
continuous representatives agreeing almost everywhere agree on
all of $[0,T]$, proving uniqueness.

In \eqref{eq:identidad-energia-galerkin-traza}, the endpoint
terms converge uniformly. Moreover,
\begin{align*}
&\left|\int_\tau^t
 \bigl(\langle\partial_\rho\mathbf{u}_N,\mathbf{u}_N\rangle_{H^{r-1}(M,\mathbf E),H^{r+1}(M,\mathbf E)}
 -\langle\partial_\rho\mathbf{u},\mathbf{u}\rangle_{H^{r-1}(M,\mathbf E),H^{r+1}(M,\mathbf E)}\bigr)
 \,d\rho\right|\\
&\quad\leq
 \|\partial_t\mathbf{u}_N-\partial_t\mathbf{u}\|_{L^2(0,T;H^{r-1}(M,\mathbf E))}
 \|\mathbf{u}_N\|_{L^2(0,T;H^{r+1}(M,\mathbf E))}\\
&\qquad+\|\partial_t\mathbf{u}\|_{L^2(0,T;H^{r-1}(M,\mathbf E))}
 \|\mathbf{u}_N-\mathbf{u}\|_{L^2(0,T;H^{r+1}(M,\mathbf E))},
\end{align*}
where the norms are taken over $(0,T)$. The right-hand side tends to
zero. Letting $N\to\infty$ gives
\eqref{eq:identidad-energia-traza-temporal} for all
$0\leq\tau\leq t\leq T$.

Finally, choose $\tau\in[0,T]$ such that
\[
 \|\mathbf{u}(\tau)\|_{H^r(M,\mathbf E)}^2
 \leq\frac{1}{T}\|\mathbf{u}\|_{L^2(0,T;H^r(M,\mathbf E))}^2
 \leq\frac{1}{T}\|\mathbf{u}\|_{L^2(0,T;H^{r+1}(M,\mathbf E))}^2.
\]
The energy identity and the same preceding argument give
\[
 \|\mathbf{u}\|_{L^\infty(0,T;H^r(M,\mathbf E))}^2
 \leq
 \left(1+\frac{1}{T}\right)
 \|\mathbf{u}\|_{L^2(0,T;H^{r+1}(M,\mathbf E))}^2
 +\|\partial_t\mathbf{u}\|_{L^2(0,T;H^{r-1}(M,\mathbf E))}^2.
\]
This proves continuity of all time evaluations with respect
to the graph norm, in particular of the trace at $t=0$.
\end{proof}

From now on, equip $\mathcal X_T^r(\mathbf{E})$ with the norm
\begin{align}
 \|\mathbf{u}\|_{\mathcal X_T^r(\mathbf{E})}^2:={}&
 \|\mathbf{u}\|_{L^\infty(0,T;H^r(M,\mathbf{E}))}^2
 +\int_0^T\|\mathbf{u}(t)\|_{H^{r+1}(M,\mathbf{E})}^2\,dt
 +\int_0^T\|\partial_t\mathbf{u}(t)\|_{H^{r-1}(M,\mathbf{E})}^2\,dt.
 \label{eq:norma-espacio-energia-parabolico}
\end{align}
The preceding estimate shows that this norm is equivalent to
\eqref{eq:norma-grafica-espacio-energia-parabolico}.

\section{Strongly parabolic linear systems}

The preceding estimates control the products arising when an
equation is differentiated. We now turn to the linear problem used at each
step of the quasilinear iteration. Set $\boldsymbol{\mathcal{H}}:=T^*M\otimes \mathbf{E}$, with its induced metric and connection.
Let
\[
 \mathbf{A}(t)\in\Gamma(\operatorname{End}\boldsymbol{\mathcal{H}}),\qquad
 \mathbf{B}(t)\in\Gamma(\operatorname{Hom}(\boldsymbol{\mathcal{H}},\mathbf{E})),\qquad
 \mathbf{C}(t)\in\Gamma(\operatorname{End}\mathbf{E})
\]
be time-dependent sections. Consider the intrinsic operator
\begin{equation}
\label{eq:operador-parabolico-lineal}
 P_t\mathbf{u}:=\nabla^*\bigl(\mathbf{A}(t)\nabla \mathbf{u}\bigr)
       +\mathbf{B}(t)(\nabla \mathbf{u})+\mathbf{C}(t)\mathbf{u}
\end{equation}
and the equation
\begin{equation}
\label{eq:problema-parabolico-lineal}
 \partial_t\mathbf{u}+P_t\mathbf{u}=\mathbf{f},
 \qquad \mathbf{u}(0)=\mathbf{u}_0.
\end{equation}
After fixing a chart and a local frame of $\mathbf{E}$, the indices
$i,j,k$ range over $\{1,\ldots,m\}$, and we write
$\displaystyle (\mathbf{A}\nabla \mathbf{u})_i=\displaystyle\sum_{j=1}^m A_i{}^j\nabla_ju$ and
$\displaystyle \mathbf{B}(\nabla \mathbf{u})=\displaystyle\sum_{j=1}^m B^j\nabla_ju$. With
$\displaystyle A^{ij}:=\displaystyle\sum_{k=1}^m g_0^{ik}A_k{}^j$, the local formula for
\eqref{eq:operador-parabolico-lineal} is
\[
 P_tu=-\sum_{i,j=1}^m\nabla_i\bigl(A^{ij}(t)\nabla_ju\bigr)
      +\sum_{j=1}^m B^j(t)\nabla_ju+C(t)u.
\]

\begin{definition}[Strong parabolicity]
\label{def:parabolicidad-fuerte-haces}
Let $(M,\mathbf{g}_0)$ be a closed Riemannian manifold and
let $\mathbf{E}\to M$ be a vector bundle with a bundle metric and compatible
connection. We say that the family of operators
\eqref{eq:operador-parabolico-lineal} is
\textbf{strongly parabolic} with constant $\lambda>0$ if $\mathbf{A}(t)$ is
self-adjoint on $\boldsymbol{\mathcal H}$ and
\begin{equation}
\label{eq:legendre-parabolicidad-fuerte}
 \left\langle \mathbf{A}(x,t)\zeta,\zeta\right\rangle_{\mathbf{g}_0,\mathbf{h}_{\mathbf{E}}}
 \geq\lambda |\zeta|_{\mathbf{g}_0,\mathbf{h}_{\mathbf{E}}}^2
\end{equation}
for every $(x,t)\in M\times[0,T]$ and
$\zeta\in T_x^*M\otimes \mathbf{E}_x$. Taking $\zeta=\xi\otimes v$ shows that this condition implies uniform
positivity of the principal Fourier symbol of the spatial operator.
When $A^{ij}=a^{ij}\operatorname{Id}_{\mathbf{E}}$, the two conditions
are equivalent.
\end{definition}

The systems we use admit a divergence-form formulation. In a
fixed chart and frame, an operator with principal part
$\displaystyle-\displaystyle\sum_{i,j=1}^m A^{ij}(\nabla^2\,\cdot\,)_{ij}$ can be written as
$\displaystyle-\displaystyle\sum_{i,j=1}^m\nabla_i(A^{ij}\nabla_j)+\displaystyle\sum_{j=1}^m\widetilde B^j\nabla_j$, with
$\displaystyle\widetilde B^j=B^j+\displaystyle\sum_{i=1}^m\nabla_iA^{ij}$ after incorporating the original term
$\displaystyle\sum_{j=1}^m B^j\nabla_j$.

For integrations by parts, use the covariant energy
\begin{equation}
\label{eq:energia-covariante-equivalente-parabolica}
 \mathscr E_r(\mathbf{u}):=\sum_{q=0}^r
 \|\nabla^q\mathbf{u}\|_{L^2(M,T^{(0,q)}(TM)\otimes \mathbf{E})}^2.
\end{equation}
There exist constants $c_r,C_r>0$ such that
\[
 c_r\|\mathbf{u}\|_{H^r(M,\mathbf{E})}^2
 \leq\mathscr E_r(\mathbf{u})
 \leq C_r\|\mathbf{u}\|_{H^r(M,\mathbf{E})}^2.
\]
The trace identity refers to the spectral norm, whereas the following
differential identities are written for $\mathscr E_r$;
integrated estimates transfer between them through this equivalence.

\begin{lemma}[Higher-order energy estimate]
\label{lem:energia-parabolica-orden-alto}
Let $(M,\mathbf{g}_0)$ be a closed Riemannian manifold and
let $\mathbf{E}\to M$ be a vector bundle with a bundle metric and compatible
connection.
Suppose that~\eqref{eq:legendre-parabolicidad-fuerte} holds and that
\begin{equation}
\label{eq:cota-coeficientes-lineales-parabolicos}
 \sup_{t\in[0,T]}
 \Bigl(
 \|\mathbf{A}(t)\|_{H^s(M,\operatorname{End}\boldsymbol{\mathcal H})}
 +\|\mathbf{B}(t)\|_{H^s(M,\operatorname{Hom}(\boldsymbol{\mathcal H},\mathbf{E}))}
 +\|\mathbf{C}(t)\|_{H^s(M,\operatorname{End}\mathbf{E})}
 \Bigr)\leq K.
\end{equation}
For $r\in\{0,1,\ldots,s\}$ there exist constants $c>0$ and $C>0$, depending only
on $r,s,\lambda,K$ and the fixed metrics and connections, such that every
smooth solution of
\eqref{eq:problema-parabolico-lineal} satisfies
\begin{equation}
\label{eq:energia-diferencial-parabolica}
 \frac{d}{dt}\mathscr E_r(\mathbf{u}(t))
 +c\|\mathbf{u}(t)\|_{H^{r+1}(M,\mathbf{E})}^2
 \leq C\mathscr E_r(\mathbf{u}(t))+C\|\mathbf{f}(t)\|_{H^{r-1}(M,\mathbf{E})}^2.
\end{equation}
Consequently,
\begin{align}
 \|\mathbf{u}\|_{L^\infty(0,T;H^r(M,\mathbf{E}))}^2
 +c\int_0^T\|\mathbf{u}(t)\|_{H^{r+1}(M,\mathbf{E})}^2\,dt
\leq
 C e^{CT}\left(
 \|\mathbf{u}_0\|_{H^r(M,\mathbf{E})}^2+
 \int_0^T\|\mathbf{f}(t)\|_{H^{r-1}(M,\mathbf{E})}^2\,dt
 \right).
 \label{eq:energia-integrada-parabolica}
\end{align}
\end{lemma}

\begin{proof}
For $q\in\{0,1,\ldots,r\}$, apply $\nabla^q$ to the equation and
pair with $\nabla^q\mathbf{u}$. The new derivative occupies the first
slot, according to Definition~\ref{derivada covariante de orden superior para haces vectoriales}.
Let $\mathbf{V}_q$ be the tensor obtained from $\nabla(\nabla^q\mathbf{u})$
by moving this first slot to the last covariant slot:
\[
 \mathbf{V}_q(X_1,\ldots,X_q,Y)
 :=(\nabla(\nabla^q\mathbf{u}))(Y,X_1,\ldots,X_q).
\]
Thus $\mathbf{A}$ acts on the pair consisting of $Y$ and the value in
$\mathbf{E}$, and $|\mathbf{V}_q|=|\nabla^{q+1}\mathbf{u}|$.
Green's formula in Theorem~\ref{teo:green-m-esima-covariante-E} gives the
principal contribution
\[
 \int_M\langle(I_{T^{(0,q)}(TM)}\otimes\mathbf{A})\mathbf{V}_q,
                    \mathbf{V}_q\rangle\,d\lambda_{\mathbf{g}_0}
 +\mathcal R_q.
\]
The integral is at least $\lambda\|\nabla^{q+1}\mathbf{u}\|_{L^2(M,T^{(0,q+1)}(TM)\otimes\mathbf E)}^2$.

Let us describe the contributions of $\mathcal R_q$. When the
$q$ derivatives are distributed between $\mathbf{A}$ and $\nabla\mathbf{u}$, the terms with
one derivative on $\mathbf{A}$ form
$[\nabla^q,\mathbf{A}](\nabla\mathbf{u})$. Reordering the remaining covariant
slots introduces commutators of two derivatives. For $q=1$ use the
curvature identity of the induced connection. If the interchanges have
been carried out up to order $q$, applying one more derivative produces two types
of terms: a derivative of a curvature coefficient already present and a
derivative of the section to which that coefficient is applied. The new
interchange contributes another curvature multiplied by a lower-order
derivative of the section. The Leibniz rule in
Proposition~\ref{prop: leibniz operador *} controls the first two types.
This is the same commutation step developed in
Lemma~\ref{lema conmutador laplaciano orden m}; here it also applies to
$T^{(0,q)}(TM)\otimes\mathbf{E}$. At order $q$, only
curvatures of the fixed structures and their derivatives up to order $q-1$ occur.
These terms therefore have bounded coefficients and contain at most
$q$ derivatives of $\mathbf{u}$ after integration by parts.

For $q=0$ there is no commutator. For $1\leq q\leq r$, estimate
\eqref{eq:moser-conmutador-asimetrico-parabolico} gives, for each
$\varepsilon>0$,
\[
 |\mathcal R_q|
 \leq C K\|\mathbf{u}\|_{H^r(M,\mathbf E)}\|\mathbf{u}\|_{H^{r+1}(M,\mathbf E)}
 \leq\varepsilon\|\mathbf{u}\|_{H^{r+1}(M,\mathbf E)}^2
       +C_\varepsilon\|\mathbf{u}\|_{H^r(M,\mathbf E)}^2.
\]
The lower-order terms obey
\[
 \|\mathbf{B}(\nabla\mathbf{u})+\mathbf{C}\mathbf{u}\|_{H^q(M,\mathbf E)}
 \leq CK\|\mathbf{u}\|_{H^{q+1}(M,\mathbf E)},\qquad q\in\{0,\ldots,r\}.
\]
Indeed, in every Leibniz product containing $|\mu|$ derivatives of
a coefficient, the available orders for applying
Lemma~\ref{lem:producto-sobolev-ordenes-distribuidos-parabolico} are
$s-|\mu|$ and $|\mu|$; their sum is $s>\displaystyle\frac{m}{2}$. This observation includes $q=0$
and avoids requiring an $L^\infty(M,\mathbf E)$ norm of the solution. After
pairing with $\nabla^q\mathbf{u}$, Young gives the same bound with
$\varepsilon$ for each of these terms.

If $r=0$, interpret the source term directly by duality:
\[
|\langle \mathbf{f},\mathbf{u}\rangle_{H^{-1}(M,\mathbf{E}),H^1(M,\mathbf{E})}|
\leq\|\mathbf{f}\|_{H^{-1}(M,\mathbf{E})}\|\mathbf{u}\|_{H^1(M,\mathbf{E})},
\]
and Young gives the required bound. Now suppose that $r\geq1$. For smooth $\mathbf{f}$,
the contribution obtained by summing the covariant energies is
\[
\mathcal B_r(\mathbf{f},\mathbf{u})
:=\sum_{q=0}^r\int_M
\langle\nabla^q\mathbf{f},\nabla^q\mathbf{u}\rangle\,d\lambda_{\mathbf{g}_0}.
\]
For $0\leq q\leq r-1$, Cauchy--Schwarz gives directly
\[
 \left|\int_M\langle\nabla^q\mathbf{f},\nabla^q\mathbf{u}\rangle
 \,d\lambda_{\mathbf{g}_0}\right|
 \leq C\|\mathbf{f}\|_{H^{r-1}(M,\mathbf E)}\|\mathbf{u}\|_{H^{r-1}(M,\mathbf E)}.
\]
In the highest-order term, write
$\nabla^r\mathbf{f}=\nabla(\nabla^{r-1}\mathbf{f})$ and apply the definition
of the covariant adjoint:
\[
 \int_M\langle\nabla^r\mathbf{f},\nabla^r\mathbf{u}\rangle
 \,d\lambda_{\mathbf{g}_0}
 =\int_M\langle\nabla^{r-1}\mathbf{f},
 \nabla^*(\nabla^r\mathbf{u})\rangle\,d\lambda_{\mathbf{g}_0}.
\]
The operator $\mathbf{u}\mapsto\nabla^*(\nabla^r\mathbf{u})$ has order
$r+1$, so its norm in
$L^2(M,T^{(0,r-1)}(TM)\otimes\mathbf E)$ is bounded by
$C\|\mathbf{u}\|_{H^{r+1}(M,\mathbf E)}$. Summing the preceding terms yields
\[
 |\mathcal B_r(\mathbf{f},\mathbf{u})|
 \leq C\|\mathbf{f}\|_{H^{r-1}(M,\mathbf{E})}
       \|\mathbf{u}\|_{H^{r+1}(M,\mathbf{E})}.
\]
For $\mathbf f\in H^{r-1}(M,\mathbf E)$ and
$\mathbf u\in H^{r+1}(M,\mathbf E)$, the expression defining the
extension is
\[
 \mathcal B_r(\mathbf f,\mathbf u)
 =\sum_{q=0}^{r-1}\int_M
   \langle\nabla_w^q\mathbf f,\nabla_w^q\mathbf u\rangle
                         \,d\lambda_{\mathbf g_0}
  +\int_M\langle\nabla_w^{r-1}\mathbf f,
                 \nabla^*(\nabla_w^r\mathbf u)\rangle
                         \,d\lambda_{\mathbf g_0}.
\]
Each pairing is between two sections in
$L^2(M,T^{(0,q)}(TM)\otimes\mathbf E)$; in the last term, take
$q=r-1$. Thus the integrals are defined.
If $\mathbf f_j\to\mathbf f$ in $H^{r-1}(M,\mathbf E)$ and
$\mathbf u_j\to\mathbf u$ in $H^{r+1}(M,\mathbf E)$ are
smooth Meyers--Serrin approximations, bilinearity and the preceding bound
give
\begin{align*}
 |\mathcal B_r(\mathbf f_j,\mathbf u_j)-\mathcal B_r(\mathbf f,\mathbf u)|
 &\leq C\|\mathbf f_j-\mathbf f\|_{H^{r-1}(M,\mathbf E)}
          \|\mathbf u_j\|_{H^{r+1}(M,\mathbf E)}\\
 &\quad+C\|\mathbf f\|_{H^{r-1}(M,\mathbf E)}
          \|\mathbf u_j-\mathbf u\|_{H^{r+1}(M,\mathbf E)}
 \longrightarrow0.
\end{align*}
Thus the extension agrees with the smooth expression, is independent
of both approximations, and satisfies the same bound.
Young then gives
\[
 |\mathcal B_r(\mathbf{f},\mathbf{u})|
\leq\varepsilon\|\mathbf{u}\|_{H^{r+1}(M,\mathbf{E})}^2
      +C\|\mathbf{f}\|_{H^{r-1}(M,\mathbf{E})}^2.
\]
Now sum all contributions. The sum of the coercive terms
controls $\|\mathbf{u}\|_{H^{r+1}(M,\mathbf E)}^2$ up to a multiple of
$\|\mathbf{u}\|_{L^2(M,\mathbf E)}^2$, which is incorporated into the right-hand side.
Choose $\varepsilon$ after this summation, smaller than the coercivity constant
divided by twice the number of terms to be absorbed.
A positive constant $c$ remains, and we obtain
\eqref{eq:energia-diferencial-parabolica}.
Integration over $[0,t]$ and Grönwall's inequality in
Lemma~\ref{lema: gronwall} prove
\eqref{eq:energia-integrada-parabolica}.
\end{proof}

\begin{theorem}[Linear existence and parabolic regularity]
\label{teo:existencia-parabolica-lineal-cerrada}
Let $(M,\mathbf{g}_0)$ be a closed Riemannian manifold and
let $\mathbf{E}\to M$ be a vector bundle with a bundle metric and compatible
connection.
Let $r\in\{0,1,\ldots,s\}$. Suppose that
\begin{align*}
 \mathbf{A}&\in C([0,T];H^s(M,\operatorname{End}\boldsymbol{\mathcal H})),\\
 \mathbf{B}&\in C([0,T];H^s(M,\operatorname{Hom}(\boldsymbol{\mathcal H},\mathbf{E}))),\\
 \mathbf{C}&\in C([0,T];H^s(M,\operatorname{End}\mathbf{E})),
\end{align*}
satisfy
\eqref{eq:legendre-parabolicidad-fuerte} with a constant independent of
$t$ and obey~\eqref{eq:cota-coeficientes-lineales-parabolicos}. For every
\[
 \mathbf{u}_0\in H^r(M,\mathbf{E}),
 \qquad \mathbf{f}\in\mathcal Y_T^r(\mathbf{E}),
\]
there exists a unique solution $\mathbf{u}\in\mathcal X_T^r(\mathbf{E})$ of
\eqref{eq:problema-parabolico-lineal}. Moreover,
\begin{equation}
\label{eq:maxima-regularidad-parabolica-lineal}
 \|\mathbf{u}\|_{\mathcal X_T^r(\mathbf{E})}
 \leq C_T\bigl(
 \|\mathbf{u}_0\|_{H^r(M,\mathbf{E})}+\|\mathbf{f}\|_{\mathcal Y_T^r(\mathbf{E})}
 \bigr),
\end{equation}
where $C_T$ can be chosen uniformly when $T$, $K$, and
$\lambda^{-1}$ remain bounded. If the coefficients, $\mathbf{f}$, and $\mathbf{u}_0$ are
smooth, then $\mathbf{u}$ is smooth on $M\times[0,T]$.
\end{theorem}

\begin{proof}
Use the same spectral basis and projection convention as in
the proof of Lemma~\ref{lem:traza-temporal-parabolica}. Order its
eigenvalues \(0\leq\mu_1\leq\mu_2\leq\cdots\), with a finite
list when applicable. Let \(V_N\) be the range of the projection
\(\Pi_N\) onto the first \(N\) vectors; if
\(d_H:=\dim L^2(M,\mathbf E)<\infty\), take
\(V_N=L^2(M,\mathbf E)\) and \(\Pi_N=I\) for \(N\geq d_H\).
Sums over \(V_N\) are restricted to the vectors in this basis.

Seek
$\displaystyle \mathbf{u}_N(t)=\displaystyle\sum_{k=1}^{\dim V_N} d_k(t)\mathbf{e}_k$ and require that, for every $\mathbf{v}\in V_N$,
\begin{align}
 (\partial_t\mathbf{u}_N,\mathbf{v})_{L^2(M,\mathbf{E})}
 &+\int_M\langle \mathbf{A}\nabla \mathbf{u}_N,\nabla \mathbf{v}\rangle_{\mathbf{g}_0,\mathbf{h}_{\mathbf{E}}}
 \,d\lambda_{\mathbf{g}_0}
 +\int_M\langle \mathbf{B}(\nabla \mathbf{u}_N)+\mathbf{C}\mathbf{u}_N,\mathbf{v}\rangle_{\mathbf{h}_{\mathbf{E}}}
 \,d\lambda_{\mathbf{g}_0}
 =\langle \mathbf{f},\mathbf{v}\rangle_{H^{-1}(M,\mathbf{E}),H^1(M,\mathbf{E})}.
 \label{eq:galerkin-parabolico}
\end{align}
For smooth data, this is a finite-dimensional linear differential
equation with continuous coefficients and therefore has a solution on all of
$[0,T]$, with $\mathbf{u}_N(0)$ equal to the orthogonal projection of $\mathbf{u}_0$ onto $V_N$.
Taking $\mathbf{v}=\mathbf{u}_N$ and using
\eqref{eq:legendre-parabolicidad-fuerte} gives, with constants
independent of $N$,
\begin{equation}
\label{eq:galerkin-energia-base}
 \|\mathbf{u}_N\|_{L^\infty(0,T;L^2(M,\mathbf{E}))}^2
 +\int_0^T\|\mathbf{u}_N(t)\|_{H^1(M,\mathbf{E})}^2\,dt
 \leq C_T\left(
 \|\mathbf{u}_0\|_{L^2(M,\mathbf{E})}^2
 +\int_0^T\|\mathbf{f}(t)\|_{H^{-1}(M,\mathbf{E})}^2\,dt
 \right).
\end{equation}
Let $\Pi_N$ be the spectral projection onto $V_N$. It is a contraction on
$H^1(M,\mathbf E)$ and $H^{-1}(M,\mathbf E)$, because it retains the first coefficients in the
spectral norms. The projected equation reads
$\partial_t\mathbf{u}_N=\Pi_N(\mathbf{f}-P_t\mathbf{u}_N)$.
The operator $P_t\colon H^1(M,\mathbf E)\to H^{-1}(M,\mathbf E)$ is uniformly bounded by
the coefficient bounds and the weak formula. The order-zero energy
then gives a uniform bound for $\partial_t\mathbf{u}_N$ in
$L^2(0,T;H^{-1}(M,\mathbf{E}))$. Reflexivity allows us to extract a subsequence
converging weakly in $L^2(0,T;H^1(M,\mathbf{E}))$ whose time derivatives
converge weakly in $L^2(0,T;H^{-1}(M,\mathbf{E}))$.

Strong convergence in $L^2(0,T;L^2(M,\mathbf{E}))$ can be seen directly. If $Q_J$
is the projection onto the first \(J\) eigensections,
in the infinite-dimensional case the \(H^1(M,\mathbf E)\) bound gives
\[
 \int_0^T\|(I-Q_J)\mathbf u_N(t)\|_{L^2(M,\mathbf E)}^2\,dt
 \leq\frac{C}{1+\mu_{J+1}}\longrightarrow0
 \quad(J\to\infty),
\]
since \(\mu_{J+1}\to\infty\). In finite dimension \(d_H\),
the tail is exactly zero for \(J\geq d_H\).

For fixed $J$, the
coordinates of $Q_J\mathbf{u}_N$ are uniformly bounded in $H^1(0,T)$ by the bound
for $\partial_t\mathbf{u}_N$ in $H^{-1}(M,\mathbf{E})$; a subsequence converges uniformly in the
finite-dimensional space $Q_JL^2(M,\mathbf{E})$. First choose $J$ large and then
a diagonal subsequence. Passing to the limit in
\eqref{eq:galerkin-parabolico} gives a weak solution $\mathbf{u}$.
Moreover, the two preceding weak convergences show that
$\mathbf{u}_N\rightharpoonup\mathbf{u}$ in the energy space
$\mathcal X_T^0(\mathbf{E})$. The trace
$\mathcal X_T^0(\mathbf{E})\to L^2(M,\mathbf{E})$ is linear and continuous by
Lemma~\ref{lem:traza-temporal-parabolica} and is therefore weakly continuous.
Since
$\mathbf{u}_N(0)=\Pi_N\mathbf{u}_0\to\mathbf{u}_0$ strongly in $L^2(M,\mathbf E)$, we
conclude that $\mathbf{u}(0)=\mathbf{u}_0$.

If the section space is finite dimensional, all \(H^r\) norms
coincide because \(\Delta_\nabla=0\). The order-zero construction
therefore already gives \(\mathcal X_T^r\) for every \(r\), without
spatial quotients. In what follows we consider the infinite-dimensional
case.

We now show that data of order $r$ give exactly the
asserted regularity. For $r=0$, the preceding construction and the equation already
give $\mathbf{u}\in\mathcal X_T^0(\mathbf{E})$. Thus suppose that
$r\geq1$. We obtain each spatial derivative using difference quotients.
Begin with the first derivative.

Let us justify the test sections depending on the solution. First fix
the spatial step $h\neq0$ and a time $T'<T$. For
$0<\varepsilon<T-T'$ use the Steklov averages
\[
 \mathbf{u}^\varepsilon(t)=\frac{1}{\varepsilon}
 \int_t^{t+\varepsilon}\mathbf{u}(\rho)\,d\rho,
 \qquad t\in[0,T'].
\]
Spatial difference quotients of $\mathbf{u}^\varepsilon$, localized in a
chart, belong to $L^2(0,T';H^1(M,\mathbf E))$ and are admissible in the averaged weak
equation. Integrate between $0$ and a time in $[0,T']$, using time cutoffs
approximating its indicator function and including the initial
term. For fixed $h$, first let $\varepsilon\to 0^{+}$.
For example, $(\mathbf{A}\nabla\mathbf{u})^\varepsilon$ and
$\mathbf{A}\nabla\mathbf{u}^\varepsilon$ both converge to
$\mathbf{A}\nabla\mathbf{u}$ strongly in
$L^2(0,T';L^2(M,T^*M\otimes\mathbf{E}))$.
This follows from continuity of time translations in $L^2$ and
the uniform bound for $\mathbf{A}$.
Fixed-step spatial difference quotients, localized in the charts,
are bounded operators on the corresponding space-time $L^2$
spaces and therefore preserve this convergence.
Moreover, $\mathbf{u}^\varepsilon(0)\to\mathbf{u}_0$ in $L^2(M,\mathbf E)$ by the order-zero trace
already obtained. This recovers the initial term for
the difference quotients. Only afterward do we let $h\to0$, using the uniform estimates
below. At the end we may take $T'\to T^{-}$.

In a chart and a trivialization of rank
$d:=\operatorname{rank}\mathbf{E}$, let $\delta_h^\ell$ be the difference quotient in the
coordinate direction $\ell$, and let $\eta$ be a smooth function compactly
supported in the chart. Choose an open subset $U'$ relatively compact in
the same chart and containing a neighborhood of $\operatorname{supp}\eta$;
the steps $h$ are chosen so that all translations used
remain in $U'$. If $d\lambda_{\mathbf{g}_0}=a(x)\,dx$ and $H(x)$ is the matrix
of the bundle metric, set $\mathsf G(x):=a(x)H(x)$. In the weak
formulation, use as test section the weighted discrete adjoint
\[
 -\mathsf G^{-1}\delta_{-h}^\ell
 \bigl(\mathsf G\eta^2\delta_h^\ell \mathbf{u}\bigr).
\]
The exact discrete identity in the chart is
\begin{equation}
\label{eq:integracion-partes-cociente-parabolico}
 \int\!\bigl\langle
 \mathsf G\delta_h^\ell v,w\bigr\rangle_{\mathbb R^d}\,dx
 =-\int\!\bigl\langle
 \mathsf Gv,
 \mathsf G^{-1}\delta_{-h}^\ell(\mathsf Gw)
 \bigr\rangle_{\mathbb R^d}\,dx.
\end{equation}
The weighted operator differs from $\delta_{-h}^\ell$ by an operator
uniformly bounded on $L^2(M,\mathbf E)$. Indeed, if
$\tau_h^\ell w(x)=w(x+h e_\ell)$, the discrete product rule is
\[
 \delta_h^\ell(fw)=(\tau_h^\ell f)\delta_h^\ell w
                         +(\delta_h^\ell f)w.
\]
In these estimates, compactly supported functions are extended by
zero to \(\mathbb R^m\). Fix a compact subset of the chart containing
all translations used, and constants
\[
 M_G:=\sup\|\mathsf G^{-1}\|,\qquad
 L_G:=\sup\|\partial_\ell\mathsf G\|.
\]
The fundamental theorem of calculus on each coordinate segment gives
\[
 \|\mathsf G(x-h e_\ell)-\mathsf G(x)\|\leq L_G|h|,
 \qquad \|\delta_{-h}^\ell\mathsf G(x)\|\leq L_G.
\]
The first difference tends to zero uniformly because \(L_G\) is fixed.
The discrete product rule gives the identity
\[
 \begin{aligned}
 \mathsf G^{-1}\delta_{-h}^\ell(\mathsf Gw)-\delta_{-h}^\ell w
 ={}&\mathsf G^{-1}(\tau_{-h}^\ell\mathsf G-\mathsf G)
                                           \delta_{-h}^\ell w\\
 &+\mathsf G^{-1}(\delta_{-h}^\ell\mathsf G)w.
 \end{aligned}
\]
Since \(\|h\delta_{-h}^\ell w\|_{L^2(\mathbb R^m)}\leq2\|w\|_{L^2(\mathbb R^m)}\),
it follows that
\[
 \|\mathsf G^{-1}\delta_{-h}^\ell(\mathsf Gw)-\delta_{-h}^\ell w\|_{L^2(\mathbb R^m)}
 \leq3M_GL_G\|w\|_{L^2(\mathbb R^m)}.
\]
Equivalence of norms on the compact set gives the same uniform bound
in \(L^2(M,\mathbf E)\).

We also quantify coercivity of the translated principal part.
In the chart, set
\[
 \mathsf H_1(x):=(g_0^{ij}(x))_{ij}\otimes\mathsf G(x),\qquad
 \mathsf H_{1,h}(x):=(g_0^{ij}(x+h e_\ell))_{ij}\otimes\mathsf G(x).
\]
The first matrix represents the metric on \(T^*M\otimes\mathbf E\)
multiplied by the density. The second appears after discrete
integration by parts: the coefficients \(g_0^{ij}\) are translated
with the flux, whereas the weight \(\mathsf G\) stays at the point
of integration.
On the fixed compact set there are \(c_1>0\) and \(L_1,C_A\geq0\) with
\[
 \mathsf H_1(x)\geq c_1I,\qquad
 \|A(x+h e_\ell,t)\|\leq C_A,
\]
\[
 \|\mathsf H_{1,h}(x)-\mathsf H_1(x+h e_\ell)\|
 +\|\mathsf H_1(x+h e_\ell)-\mathsf H_1(x)\|\leq L_1|h|.
\]
The last bound follows from the fundamental theorem of calculus
applied to the metric coefficients and their products.
Here \(C_A\) is controlled by \(K\) and the fixed structures.
For the column \(Z\) of covariant derivatives, parabolicity at
\(x+h e_\ell\) implies
\[
 \begin{aligned}
 Z^{\mathsf T}\mathsf H_{1,h}(x)A(x+h e_\ell,t)Z
 &\geq\lambda Z^{\mathsf T}\mathsf H_1(x)Z
                -(\lambda+C_A)L_1|h|\,|Z|^2\\
 &\geq(\lambda-C_0|h|)Z^{\mathsf T}\mathsf H_1(x)Z,
 \qquad C_0:=\frac{(\lambda+C_A)L_1}{c_1}.
 \end{aligned}
\]
Connection terms outside this principal form are included among the
lower-order terms. Choose a step satisfying
\[
 0<|h|<\min\left\{h_0,\frac{\lambda}{2(1+C_0)}\right\},
\]
where \(h_0>0\) keeps the segments inside the chart. Then
\(\lambda-C_0|h|\geq\displaystyle\frac{\lambda}{2}\), and the
principal integral controls
\[
 \frac{\lambda}{2}\int_M\eta^2
 |\nabla\delta_h^\ell\mathbf u|_{\mathbf g_0,\mathbf h_{\mathbf E}}^2
                         \,d\lambda_{\mathbf g_0}.
\]
All these constants are independent of the step \(h\).

The time term is
\[
 \frac{1}{2}\frac{d}{dt}
 \int_M\eta^2|\delta_h^\ell \mathbf{u}|_{\mathbf{h}_{\mathbf{E}}}^2\,d\lambda_{\mathbf{g}_0}.
\]
Terms in which the difference quotient falls on the local representation of $\mathbf{A}$,
on $\eta$, or on the
connection coefficients are bounded by
\begin{align*}
 \frac{\lambda}{4}
 \int_M\eta^2|\nabla\delta_h^\ell \mathbf{u}|_{\mathbf{g}_0,\mathbf{h}_{\mathbf{E}}}^2\,d\lambda_{\mathbf{g}_0}+C\int_{U'}
\bigl(|\nabla \mathbf{u}|_{\mathbf{g}_0,\mathbf{h}_{\mathbf{E}}}^2+|\mathbf{u}|_{\mathbf{h}_{\mathbf{E}}}^2+|\mathbf{f}|_{\mathbf{h}_{\mathbf{E}}}^2\bigr)
 \,d\lambda_{\mathbf{g}_0}.
\end{align*}
The bound is uniform for $|h|$ smaller than the distance between
$\operatorname{supp}\eta$ and the edge of the chart. Integrating over $[0,T]$
and letting $h\to0$ gives one additional spatial derivative in a smaller
chart.

To iterate, we do not assume in advance that $\nabla \mathbf{u}$ is bounded. Apply
difference quotients of order $q$, $1\leq q\leq r$, and use as test sections
the backward quotients of order $q$ of the same localized expression. The
commutator with the principal part follows from the discrete product
rule. For $\gamma\in\mathbb N_0^m$ and a common step in each direction,
writing $\tau_h^\gamma$ for translation by $h\gamma$, we have
\[
 \delta_h^\gamma(fw)=
 \sum_{\substack{\mu\in\mathbb N_0^m\\\mu\leq\gamma}}
 \binom\gamma\mu
 (\tau_h^{\gamma-\mu}\delta_h^\mu f)
 (\delta_h^{\gamma-\mu}w).
\]
This identity follows by applying the first-order rule once for
each direction in the multiindex. Translations are bounded in the localized
norms, and quotients of order $|\mu|$ are controlled by derivatives
of that order in $U'$. Thus the commutator terms satisfy the same
bounds as the differential expressions to which they converge:
\begin{equation}
\label{eq:conmutadores-regularidad-parabolica}
 (\partial^\mu A)\,
 \partial^{\gamma-\mu+1}\mathbf{u},
 \qquad \mu,\gamma\in\mathbb N_0^m,\quad 0<\mu\leq\gamma,\quad |\gamma|=q,
\end{equation}
and the connection terms have strictly lower order. The asymmetric
estimate~\eqref{eq:moser-conmutador-asimetrico-parabolico} shows that the sum of
all terms~\eqref{eq:conmutadores-regularidad-parabolica} is
bounded by
$C\|\mathbf{A}\|_{H^s(M,\operatorname{End}\boldsymbol{\mathcal H})}
\|\mathbf{u}\|_{H^{|\gamma|}(M,\mathbf{E})}$. The first-order and zeroth-order terms
obey
\[
 \|\mathbf{B}(\nabla \mathbf{u})+\mathbf{C}\mathbf{u}\|_{H^{q-1}(M,\mathbf{E})}
 \leq C K\|\mathbf{u}\|_{H^q(M,\mathbf{E})},
\]
while the sum of the source terms is the preceding covariant form
$\mathcal B_q(\mathbf{f},\mathbf{u})$. Moving one derivative in the highest-order term
and then using spectral duality
\eqref{eq:emparejamiento-terna-parabolica} gives
\[
|\mathcal B_q(\mathbf{f},\mathbf{u})|
\leq\frac{c}{4}\|\mathbf{u}\|_{H^{q+1}(M,\mathbf{E})}^2
+C\|\mathbf{f}\|_{H^{q-1}(M,\mathbf{E})}^2.
\]
Let us specify the induction that passes from these local estimates to
global regularity. For \(0\leq q\leq r\), let \(\mathcal P_q\) be the
assertion
\[
\mathbf{u}\in C([0,T];H^q(M,\mathbf{E}))
\cap L^2(0,T;H^{q+1}(M,\mathbf{E}))
\]
together with the energy estimate of order \(q\) written below.
The assertion \(\mathcal P_0\) is
\eqref{eq:galerkin-energia-base} after passage to the limit. The difference-quotient
calculation just carried out proves \(\mathcal P_1\): the uniform bound
for \(\delta_h^\ell\mathbf{u}\) allows a weak limit to be extracted,
and the difference-quotient characterization identifies that limit
with \(\partial_\ell\mathbf{u}\).

Now suppose that
\(\mathcal P_0,\dots,\mathcal P_{q-1}\) have been proved, where \(2\leq q\leq r\).
Apply a difference quotient of order \(q\) to the equation and test against its
localized backward quotient. When all quotients fall on
\(\mathbf{u}\), parabolicity gives the coercive term of order \(q+1\).
If at least one falls on \(\mathbf{A}\), we obtain a difference-quotient and translation version of one of the
terms in \eqref{eq:conmutadores-regularidad-parabolica}; the preceding Moser
estimate and Young allow the part containing
\(\|\mathbf{u}\|_{H^{q+1}(M,\mathbf E)}^2\) to be absorbed. Quotients falling on
\(\mathbf{B},\mathbf{C}\), connection coefficients, or cutoffs
contain at most \(q\) derivatives of \(\mathbf{u}\) and are controlled
by \(\mathcal P_{q-1}\) and the uniform coefficient bound. The
source term is controlled by \(\mathcal B_q\), as shown
above. Finally, quotients of the initial data are uniformly
bounded by \(\|\mathbf{u}_0\|_{H^q(M,\mathbf E)}\).

Choose a finite cover by charts and a subordinate partition of
unity. Summing the local estimates, terms in which
derivatives reach the cutoffs have lower order and are already covered
by \(\mathcal P_{q-1}\). The resulting bounds are uniform in the quotient
step; letting \(h\to0\), weak lower semicontinuity of the norm
gives
\(\mathbf{u}\in L^2(0,T;H^{q+1}(M,\mathbf{E}))\). The equation itself,
viewed as an identity in \(H^{q-1}(M,\mathbf{E})\), then gives
\(\partial_t\mathbf{u}\in L^2(0,T;H^{q-1}(M,\mathbf{E}))\). Applying
Lemma~\ref{lem:traza-temporal-parabolica} to the Hilbert triple
\[
H^{q+1}(M,\mathbf{E})\lhook\joinrel\longrightarrow
H^q(M,\mathbf{E})\lhook\joinrel\longrightarrow
H^{q-1}(M,\mathbf{E})
\]
yields \(\mathbf{u}\in C([0,T];H^q(M,\mathbf{E}))\). This proves
\(\mathcal P_q\). The inductive step is complete, and for each
\(1\leq q\leq r\) we obtain
\begin{equation}
\label{eq:induccion-cocientes-parabolica}
 \frac{d}{dt}\mathscr E_q(\mathbf{u}(t))
 +c\|\mathbf{u}(t)\|_{H^{q+1}(M,\mathbf{E})}^2
 \leq C\mathscr E_q(\mathbf{u}(t))+C\|\mathbf{f}(t)\|_{H^{q-1}(M,\mathbf{E})}^2.
\end{equation}
Each step uses only the coefficient bound and regularity from
the preceding step. Integrating the last inequality gives
\begin{align}
 \|\mathbf{u}\|_{L^\infty(0,T;H^r(M,\mathbf{E}))}^2
 +\int_0^T\|\mathbf{u}(t)\|_{H^{r+1}(M,\mathbf{E})}^2\,dt
 \leq C_T\left(
 \|\mathbf{u}_0\|_{H^r(M,\mathbf{E})}^2
 +\int_0^T\|\mathbf{f}(t)\|_{H^{r-1}(M,\mathbf{E})}^2\,dt
 \right).
 \label{eq:regularidad-cocientes-parabolica-global}
\end{align}
Control up to $t=0$ follows by applying the same calculation to
approximations with smooth initial data and using the term
$\|\mathbf{u}_0\|_{H^r(M,\mathbf{E})}$ when integrating the energy identity. No additional compatibility condition arises, since $M$ is closed.

The equation and continuity of multiplication imply
\begin{align*}
 \|\partial_t\mathbf{u}\|_{L^2(0,T;H^{r-1}(M,\mathbf{E}))}
 \leq C\bigl(&\|\mathbf{u}\|_{L^2(0,T;H^{r+1}(M,\mathbf{E}))}
 +\|\mathbf{u}\|_{L^2(0,T;H^r(M,\mathbf{E}))}
 +\|\mathbf{f}\|_{L^2(0,T;H^{r-1}(M,\mathbf{E}))}\bigr).
\end{align*}
This proves~\eqref{eq:maxima-regularidad-parabolica-lineal} for smooth data and
coefficients.

We turn to coefficients and data with the regularity in the statement. Choose
smooth approximations $\mathbf{A}_j,\mathbf{B}_j,\mathbf{C}_j,\mathbf{f}_j$ and
$\mathbf{u}_{0,j}$ converging in the stated spaces. Replace
$\mathbf{A}_j$ by its self-adjoint part. Since
$H^s(M)\hookrightarrow C^1(M)$, convergence of $\mathbf{A}_j$ is uniform;
thus, from some index onward,
\[
 \langle\mathbf{A}_j(x,t)\zeta,\zeta\rangle
 \geq\frac{\lambda}{2}|\zeta|^2
\]
for every $(x,t)$ and every $\zeta\in T_x^*M\otimes\mathbf{E}_x$. The smooth
solutions $\mathbf{u}_j$ satisfy a uniform bound in
$\mathcal X_T^r(\mathbf{E})$ by
\eqref{eq:energia-integrada-parabolica} and the equation itself.

Reflexivity of the Bochner spaces allows us to extract a subsequence
such that
\begin{align*}
 \mathbf{u}_j&\rightharpoonup\mathbf{u}
 &&\text{in }L^2(0,T;H^{r+1}(M,\mathbf{E})),\\
 \partial_t\mathbf{u}_j&\rightharpoonup\partial_t\mathbf{u}
 &&\text{in }L^2(0,T;H^{r-1}(M,\mathbf{E})).
\end{align*}
To pass to the limit in the principal term, write
\[
 \mathbf{A}_j\nabla\mathbf{u}_j-\mathbf{A}\nabla\mathbf{u}
 =(\mathbf{A}_j-\mathbf{A})\nabla\mathbf{u}_j
 +\mathbf{A}\nabla(\mathbf{u}_j-\mathbf{u}).
\]
The first part converges strongly to zero in
$L^2(0,T;H^{r-1}(M,\mathbf E))$ by the product estimates, and the second converges
weakly. The terms containing $\mathbf{B}_j$ and $\mathbf{C}_j$ are handled in the
same way. Testing the equation against time-dependent smooth sections
therefore gives the limiting problem.

The trace map
$\mathcal X_T^r(\mathbf{E})\to H^r(M,\mathbf{E})$ is linear and continuous by
Lemma~\ref{lem:traza-temporal-parabolica} and is therefore weakly
continuous. Since $\mathbf{u}_j(0)=\mathbf{u}_{0,j}\to\mathbf{u}_0$, we conclude
that $\mathbf{u}(0)=\mathbf{u}_0$. Weak lower semicontinuity preserves
estimate \eqref{eq:maxima-regularidad-parabolica-lineal}. If two solutions
have the same data, the order-zero energy estimate applied to
their difference gives uniqueness. Finally, when all data are smooth,
repeat the estimate at every spatial order; the equation expresses
time derivatives in terms of the spatial derivatives already
controlled and gives joint smoothness.
\end{proof}

\begin{proposition}[Stability with respect to coefficients]
\label{prop:estabilidad-coeficientes-parabolicos}
Let $(M,\mathbf{g}_0)$ be a closed Riemannian manifold and
let $\mathbf{E}\to M$ be a vector bundle with a bundle metric and compatible
connection.
Let $P_t$ and $\widetilde P_t$ be two families satisfying the hypotheses of
Theorem~\ref{teo:existencia-parabolica-lineal-cerrada} with the same bounds
$K$ and $\lambda$. Let $\mathbf{u},\widetilde{\mathbf{u}}\in\mathcal X_T^s(\mathbf{E})$ be the solutions
corresponding to $(\mathbf{u}_0,\mathbf{f})$ and
$(\widetilde{\mathbf{u}}_0,\widetilde{\mathbf{f}})$. Then
\begin{equation}
\begin{aligned}
 &\|\mathbf{u}-\widetilde{\mathbf{u}}\|_{\mathcal X_T^{s-1}(\mathbf{E})}\\
 &\quad\leq C_T\biggl(
 \|\mathbf{u}_0-\widetilde{\mathbf{u}}_0\|_{H^{s-1}(M,\mathbf{E})}
 +\|\mathbf{f}-\widetilde{\mathbf{f}}\|_{\mathcal Y_T^{s-1}(\mathbf{E})}\\
 &\qquad+
 \|\mathbf{A}-\widetilde{\mathbf{A}}\|_{C([0,T];H^{s-1}(M,\operatorname{End}\boldsymbol{\mathcal H}))}
 \|\widetilde{\mathbf{u}}\|_{L^2(0,T;H^{s+1}(M,\mathbf{E}))}\\
 &\qquad+
 \|\mathbf{B}-\widetilde{\mathbf{B}}\|_{C([0,T];H^{s-1}(M,\operatorname{Hom}(\boldsymbol{\mathcal H},\mathbf{E})))}
 \|\widetilde{\mathbf{u}}\|_{L^2(0,T;H^s(M,\mathbf{E}))}\\
 &\qquad+
 \|\mathbf{C}-\widetilde{\mathbf{C}}\|_{C([0,T];H^{s-1}(M,\operatorname{End}\mathbf{E}))}
 \|\widetilde{\mathbf{u}}\|_{L^2(0,T;H^s(M,\mathbf{E}))}\biggr).
\end{aligned}
\label{eq:estabilidad-coeficientes-parabolicos}
\end{equation}
\end{proposition}

\begin{proof}
Set $\mathbf{w}=\mathbf{u}-\widetilde{\mathbf{u}}$. Subtracting the equations and retaining $P_t$ on
the left-hand side gives
\begin{align*}
 \partial_t\mathbf{w}+P_t\mathbf{w}
 =\mathbf{f}-\widetilde{\mathbf{f}}
 -\nabla^*\bigl((\mathbf{A}-\widetilde{\mathbf{A}})\nabla\widetilde{\mathbf{u}}\bigr)-(\mathbf{B}-\widetilde{\mathbf{B}})(\nabla\widetilde{\mathbf{u}})
 -(\mathbf{C}-\widetilde{\mathbf{C}})\widetilde{\mathbf{u}}.
\end{align*}
Since $s-1>\displaystyle\frac{m}{2}$, product estimate
\eqref{eq:moser-producto-parabolico} shows that the right-hand side belongs
to $L^2(0,T;H^{s-2}(M,\mathbf{E}))$ and is bounded by the terms appearing in
\eqref{eq:estabilidad-coeficientes-parabolicos}. Apply
Theorem~\ref{teo:existencia-parabolica-lineal-cerrada} with $r=s-1$.
\end{proof}

\section{Quasilinear equations}

The linear theory can be applied with coefficients frozen at a given
section. The nonlinear solution will be constructed by repeating this procedure and
showing that, on a sufficiently short interval, the iteration stays
in a uniformly parabolic region and becomes contractive in a weaker
norm. Let $\mathcal U\subseteq \mathbf{E}$ be an open subset of a bundle and consider two smooth bundle
maps
\begin{align*}
 A&\colon\mathcal U\longrightarrow
 \operatorname{End}(T^*M\otimes \mathbf{E}),\\
 F&\colon J^1(\mathcal U)\longrightarrow \mathbf{E}.
\end{align*}
Here $J^1(\mathcal U)$ denotes the open subset of the bundle
$\mathbf{E}\oplus(T^*M\otimes \mathbf{E})$ consisting of pairs $(z,p)$ with $z\in\mathcal U$.
We study
\begin{equation}
\label{eq:ecuacion-parabolica-cuaslineal}
 \partial_t\mathbf{u}+
 \nabla^*\bigl(A(\mathbf{u})\nabla \mathbf{u}\bigr)
 =F(j^1\mathbf{u}),
 \qquad \mathbf{u}(0)=\mathbf{u}_0.
\end{equation}
For $\xi\in T_x^*M$, let
$\iota_\xi\colon\mathbf{E}_x\to T_x^*M\otimes\mathbf{E}_x$ be the map
$v\mapsto\xi\otimes v$. With $\mathbf{u}$ fixed, the principal Fourier
symbol of the spatial operator
$\mathbf{v}\mapsto\nabla^*(A(\mathbf{u})\nabla\mathbf{v})$ is
$\iota_\xi^*A(\mathbf{u}(x))\iota_\xi$; its differential symbol,
with the convention of
Definition~\ref{def:simbolo-principal-pdo-haces}, is the negative of
this expression, by~\eqref{eq:conversion-simbolos-diferencial-pseudo}. In a chart and local frame, the equation
takes the form
\[
 \partial_tu-\sum_{i,j=1}^m\nabla_i\bigl(A^{ij}(x,u)\nabla_ju\bigr)
 =F(x,u,\nabla u).
\]
Expanding the divergence, the horizontal derivative $\nabla_xA$ and
vertical derivative $D_zA(u)[\nabla u]$ contribute only lower-order
terms.

\begin{theorem}[Short-time quasilinear existence]
\label{teo:existencia-parabolica-cuaslineal}
Let $(M,\mathbf{g}_0)$ be a closed Riemannian manifold, let
$\mathbf{E}\to M$ be a vector bundle with a bundle metric and compatible
connection, and let $\mathcal U\subseteq\mathbf{E}$ be an open subset of the bundle.
Let $\mathbf{u}_0\in H^s(M,\mathbf{E})$ be a section whose image is contained in
$\mathcal U$. Suppose that there exists $\lambda_0>0$ such that
\begin{equation}
\label{eq:parabolicidad-inicial-cuaslineal}
 \left\langle A(\mathbf{u}_0(x))\zeta,\zeta\right\rangle_{\mathbf{g}_0,\mathbf{h}_{\mathbf{E}}}
 \geq\lambda_0|\zeta|_{\mathbf{g}_0,\mathbf{h}_{\mathbf{E}}}^2
\end{equation}
for every $x\in M$ and $\zeta\in T_x^*M\otimes \mathbf{E}_x$, and that
$A(z)$ is self-adjoint for every $z\in\mathcal U$. Then there exists $T>0$ for which
\eqref{eq:ecuacion-parabolica-cuaslineal} has a unique solution
\[
 \mathbf{u}\in\mathcal X_T^s(\mathbf{E})
 \quad\hbox{with}\quad
 \mathbf{u}(x,t)\in\mathcal U.
\]
The time $T$ can be chosen uniformly as $\mathbf{u}_0$ ranges over a
bounded subset of $H^s(M,\mathbf{E})$ whose image stays a positive distance from
the complement of $\mathcal U$ and for which the constant in
\eqref{eq:parabolicidad-inicial-cuaslineal} has a positive lower bound.
\end{theorem}

\begin{proof}
The image of $\mathbf u_0$ is compact because $s>\displaystyle\frac{m}{2}$ and $M$ is compact.
For $b>0$, define the closed fiberwise tube
\[
 T_b(\mathbf u_0):=
 \{z\in\mathbf E_x\mid x\in M,\
                     |z-\mathbf u_0(x)|_{\mathbf h_{\mathbf E}}\leq b\}.
\]
It is compact, being the image of the closed disk bundle of radius $b$
over the compact base under continuous translation by $\mathbf u_0$.
We can choose $\varepsilon>0$ such that
\[
 T_{2\varepsilon}(\mathbf u_0)\subseteq\mathcal U,
 \qquad
 \sup_{x\in M}\ \sup_{\substack{z\in\mathbf E_x\\
          |z-\mathbf u_0(x)|\leq2\varepsilon}}
 \|A(z)-A(\mathbf u_0(x))\|_{\mathrm{op}}<\frac{\lambda_0}{2}.
\]
Here the operator norm is induced by $\mathbf g_0$ and
$\mathbf h_{\mathbf E}$ on $T_x^*M\otimes\mathbf E_x$.
To justify the choice, if it failed for radii tending to zero, we could
choose $x_j\in M$ and $z_j\in\mathbf E_{x_j}$ with
$|z_j-\mathbf u_0(x_j)|\to0$ violating one of the conditions. A subsequence
converges in $T_1(\mathbf u_0)$ to $\mathbf u_0(x)$.
Openness of $\mathcal U$ and continuity of $A$ contradict, respectively,
each of these violations.
Set
\[
 \mathcal V:=\{z\in\mathbf E_x\mid x\in M,\
                   |z-\mathbf u_0(x)|_{\mathbf h_{\mathbf E}}<\varepsilon\}.
\]
This is open, contains $\mathbf u_0(M)$, and has compact closure contained
in $\mathcal U$. For $z\in\overline{\mathcal V}\cap\mathbf E_x$, the operator
bound and \eqref{eq:parabolicidad-inicial-cuaslineal} give
\begin{equation}
\label{eq:parabolicidad-vecindad-cuaslineal}
 \left\langle A(z)\zeta,\zeta\right\rangle_{\mathbf{g}_0,\mathbf{h}_{\mathbf{E}}}
 \geq\frac{\lambda_0}{2}|\zeta|_{\mathbf{g}_0,\mathbf{h}_{\mathbf{E}}}^2.
\end{equation}

Define $\mathbf{u}^{(0)}(t)=\mathbf{u}_0$ and, once $\mathbf{u}^{(k)}$ is constructed, solve the
linear problem
\begin{equation}
\label{eq:iteracion-parabolica-cuaslineal}
 \partial_t\mathbf{u}^{(k+1)}+
 \nabla^*\bigl(A(\mathbf{u}^{(k)})\nabla \mathbf{u}^{(k+1)}\bigr)
 =F(j^1\mathbf{u}^{(k)}),
 \qquad \mathbf{u}^{(k+1)}(0)=\mathbf{u}_0.
\end{equation}
We show simultaneously that the iteration is well defined and
contractive.

Let $D_s=c_s^{-1}$, with $c_s$ as in the equivalence following
\eqref{eq:energia-covariante-equivalente-parabolica}. Fix $R$ large enough
that $R^2>2D_s\mathscr E_s(\mathbf{u}_0)+1$. Suppose inductively that
\begin{equation}
\label{eq:bola-iteracion-cuaslineal}
 \|\mathbf{u}^{(k)}\|_{L^\infty(0,T;H^s(M,\mathbf{E}))}\leq R,
 \qquad \mathbf{u}^{(k)}(M\times[0,T])\subseteq\mathcal V.
\end{equation}
Lemma~\ref{lem:moser-parabolico-haces} gives bounds independent of $k$
for $A(\mathbf{u}^{(k)})$ in
$C([0,T];H^s(M,\operatorname{End}\boldsymbol{\mathcal{H}}))$ and for
$F(j^1\mathbf{u}^{(k)})$ in $C([0,T];H^{s-1}(M,\mathbf{E}))$. Together with
\eqref{eq:parabolicidad-vecindad-cuaslineal},
Theorem~\ref{teo:existencia-parabolica-lineal-cerrada} gives
$\mathbf{u}^{(k+1)}\in\mathcal X_T^s(\mathbf{E})$. Differential inequality
\eqref{eq:energia-diferencial-parabolica}, before replacement by the coarser
bound~\eqref{eq:maxima-regularidad-parabolica-lineal}, gives a constant
$C_R$ independent of $k$ and $T\leq1$ such that
\begin{align}
 \|\mathbf{u}^{(k+1)}\|_{L^\infty(0,T;H^s(M,\mathbf{E}))}^2
 &\leq D_s e^{C_RT}\bigl(\mathscr E_s(\mathbf{u}_0)+C_RT\bigr),
 \label{eq:cota-sup-iteracion-cuaslineal}\\
 \|\mathbf{u}^{(k+1)}\|_{\mathcal X_T^s(\mathbf{E})}&\leq C_R.
 \label{eq:cota-uniforme-iteracion-cuaslineal}
\end{align}
The first bound tends to $D_s\mathscr E_s(\mathbf{u}_0)$ as $T\to0$. Thus,
after reducing $T$, its right-hand side is less than $R^2$.

Moreover, by Lemma~\ref{lem:traza-temporal-parabolica},
\begin{align}
 \|\mathbf{u}^{(k+1)}(t)-\mathbf{u}_0\|_{H^{s-1}(M,\mathbf{E})}
 &\leq\int_0^t
 \|\partial_\rho \mathbf{u}^{(k+1)}(\rho)\|_{H^{s-1}(M,\mathbf{E})}\,d\rho
 \nonumber\\
 &\leq T^{\frac{1}{2}}
 \|\partial_t\mathbf{u}^{(k+1)}\|_{L^2(0,T;H^{s-1}(M,\mathbf{E}))}.
 \label{eq:cercania-inicial-iteracion-cuaslineal}
\end{align}
Let $C_{\mathrm{emb}}>0$ be an embedding constant for
$H^{s-1}(M,\mathbf E)\hookrightarrow C^0(M,\mathbf E)$. Increase $C_R$
if necessary so that $C_R\geq1$, and reduce $T$ so that, in addition to
the preceding bound by $R^2$,
\[
 C_{\mathrm{emb}}C_R\sqrt T<\frac{\varepsilon}{2}.
\]
Estimates \eqref{eq:cota-uniforme-iteracion-cuaslineal} and
\eqref{eq:cercania-inicial-iteracion-cuaslineal} then give
\[
 \sup_{0\leq t\leq T}
 \|\mathbf u^{(k+1)}(t)-\mathbf u_0\|_{C^0(M,\mathbf E)}
 \leq C_{\mathrm{emb}}C_R\sqrt T<\frac{\varepsilon}{2}.
\]
The image of $\mathbf u^{(k+1)}$ remains in the tube
$T_{\displaystyle\frac{\varepsilon}{2}}(\mathbf u_0)\subseteq\mathcal V$. The same $T$ works
for every index $k$, closing the induction.

Set $\mathbf{w}^{(k+1)}=\mathbf{u}^{(k+1)}-\mathbf{u}^{(k)}$ and
$\mathbf{A}_k=A(\mathbf{u}^{(k)})$, $\mathbf{F}_k=F(j^1\mathbf{u}^{(k)})$. The difference of two
consecutive equations is
\begin{align}
 \partial_t\mathbf{w}^{(k+1)}+
 \nabla^*(\mathbf{A}_k\nabla \mathbf{w}^{(k+1)})
 ={}&-\nabla^*\bigl((\mathbf{A}_k-\mathbf{A}_{k-1})\nabla \mathbf{u}^{(k)}\bigr)
 \nonumber\\
 &+\mathbf{F}_k-\mathbf{F}_{k-1},
 \qquad \mathbf{w}^{(k+1)}(0)=0.
 \label{eq:diferencia-iteracion-cuaslineal}
\end{align}
Contraction is obtained in the lower-order energy space
\[
 \mathcal Z_T(\mathbf{E}):=C([0,T];L^2(M,\mathbf{E}))\cap L^2(0,T;H^1(M,\mathbf{E})),
\]
with
\begin{equation}
\label{eq:norma-Z-iteracion-cuaslineal}
 \|\mathbf{w}\|_{\mathcal Z_T(\mathbf{E})}^2
 :=\|\mathbf{w}\|_{L^\infty(0,T;L^2(M,\mathbf{E}))}^2
   +\int_0^T\|\mathbf{w}(t)\|_{H^1(M,\mathbf{E})}^2\,dt.
\end{equation}
The preceding bounds and the embedding $H^s(M,\mathbf{E})\hookrightarrow C^1(M,\mathbf{E})$
confine the first jets of the iteration to a compact subset of
$J^1(\mathcal U)$. Using cutoff functions as in the proof of
Lemma~\ref{lem:moser-parabolico-haces}, extend $A$ and $F$ to smooth maps
defined on the full vector fibers of a neighborhood of this
compact set. The mean value theorem applied to these extensions allows us to
write
\begin{align}
 \mathbf{A}_k-\mathbf{A}_{k-1}&=\boldsymbol{\mathcal A}_k[\mathbf{w}^{(k)}],
 \label{eq:diferencia-A-iteracion-L2}\\
 \mathbf{F}_k-\mathbf{F}_{k-1}&=\boldsymbol{\mathcal B}_k[\mathbf{w}^{(k)}]
              +\boldsymbol{\mathcal C}_k[\nabla \mathbf{w}^{(k)}].
 \label{eq:diferencia-F-iteracion-L2}
\end{align}
For every $t\in[0,T]$, these coefficients belong to
\begin{align*}
 \boldsymbol{\mathcal A}_k(t)&\in
 H^s\bigl(M,\operatorname{Hom}(\mathbf{E},\operatorname{End}\boldsymbol{\mathcal H})\bigr),\\
 \boldsymbol{\mathcal B}_k(t)&\in H^{s-1}(M,\operatorname{End}\mathbf{E}),\\
 \boldsymbol{\mathcal C}_k(t)&\in
 H^{s-1}(M,\operatorname{Hom}(\boldsymbol{\mathcal H},\mathbf{E})),
\end{align*}
and depend continuously on $t$ in these spaces.
The bounds~\eqref{eq:bola-iteracion-cuaslineal} and the embedding
$H^s(M,\mathbf{E})\hookrightarrow C^2(M,\mathbf{E})$ give uniform bounds for the pointwise norms of
$\boldsymbol{\mathcal A}_k,\boldsymbol{\mathcal B}_k,\boldsymbol{\mathcal C}_k$,
$\nabla\boldsymbol{\mathcal C}_k$, and
$\nabla \mathbf{u}^{(k)}$.
Here $\nabla\boldsymbol{\mathcal C}_k$ is a section of
$T^*M\otimes\operatorname{Hom}(\boldsymbol{\mathcal{H}},\mathbf{E})$ with the induced connection.

Take the $L^2(M,\mathbf{E})$ inner product of
\eqref{eq:diferencia-iteracion-cuaslineal} with $\mathbf{w}^{(k+1)}$. Integrate by
parts the term containing $\mathbf{A}_k-\mathbf{A}_{k-1}$ and also the term with
$\nabla \mathbf{w}^{(k)}$. Since $\boldsymbol{\mathcal C}_k$ and
$\nabla\boldsymbol{\mathcal C}_k$ are
uniformly bounded, the latter satisfies
\begin{align*}
 \left|\int_M\left\langle
 \boldsymbol{\mathcal C}_k[\nabla \mathbf{w}^{(k)}],\mathbf{w}^{(k+1)}
 \right\rangle_{\mathbf{h}_{\mathbf{E}}}\,d\lambda_{\mathbf{g}_0}\right|
 \leq C_R\|\mathbf{w}^{(k)}\|_{L^2(M,\mathbf{E})}
 \bigl(\|\nabla \mathbf{w}^{(k+1)}\|_{L^2(M,\boldsymbol{\mathcal{H}})}+\|\mathbf{w}^{(k+1)}\|_{L^2(M,\mathbf{E})}\bigr).
\end{align*}
Parabolicity and Young's inequality give
\begin{align}
 \frac{d}{dt}\|\mathbf{w}^{(k+1)}(t)\|_{L^2(M,\mathbf{E})}^2
 +\frac{\lambda_0}{2}
 \|\nabla \mathbf{w}^{(k+1)}(t)\|_{L^2(M,\boldsymbol{\mathcal{H}})}^2
 \leq C_R\|\mathbf{w}^{(k+1)}(t)\|_{L^2(M,\mathbf{E})}^2+C_R\|\mathbf{w}^{(k)}(t)\|_{L^2(M,\mathbf{E})}^2.
 \label{eq:energia-contraccion-L2-cuaslineal}
\end{align}
Since $\mathbf{w}^{(k+1)}(0)=0$, Grönwall's lemma~\ref{lema: gronwall} implies
\begin{equation}
\label{eq:contraccion-iteracion-cuaslineal}
 \|\mathbf{w}^{(k+1)}\|_{\mathcal Z_T(\mathbf{E})}^2
 \leq C_RTe^{C_RT}
       \|\mathbf{w}^{(k)}\|_{C([0,T];L^2(M,\mathbf{E}))}^2
 \leq C_RTe^{C_RT}\|\mathbf{w}^{(k)}\|_{\mathcal Z_T(\mathbf{E})}^2.
\end{equation}
Reduce $T$ until $q:=C_RTe^{C_RT}<1$. Then
\[
 \|\mathbf{w}^{(k+1)}\|_{\mathcal Z_T(\mathbf{E})}
 \leq q^{\frac{1}{2}}\|\mathbf{w}^{(k)}\|_{\mathcal Z_T(\mathbf{E})},
\]
For $k=1$ the bound
$\|\mathbf{w}^{(k)}\|_{\mathcal Z_T}\leq
q^{\displaystyle\frac{k-1}{2}}\|\mathbf{w}^{(1)}\|_{\mathcal Z_T}$ is an equality.
If it holds for $k$, the preceding inequality gives it for $k+1$ after multiplication
by $q^{\displaystyle\frac{1}{2}}$. Thus, summing the geometric series,
\[
 \sum_{k=1}^{\infty}\|\mathbf{w}^{(k)}\|_{\mathcal Z_T(\mathbf{E})}
 \leq\frac{\|\mathbf{w}^{(1)}\|_{\mathcal Z_T(\mathbf{E})}}{1-q^{\frac{1}{2}}}<\infty.
\]
The space $\mathcal Z_T(\mathbf{E})$ is complete as the intersection of the two
Banach spaces defining it, with the sum norm. Therefore,
$(\mathbf{u}^{(k)})$ is Cauchy in $\mathcal Z_T(\mathbf{E})$.

The uniform bounds for $(\mathbf{u}^{(k)})_{k\geq1}$ in $\mathcal X_T^s(\mathbf{E})$ improve this
convergence. The section $\mathbf{u}^{(0)}(t)=\mathbf{u}_0$ is used only to start the
iteration and need not belong to $L^2(0,T;H^{s+1}(M,\mathbf{E}))$. At each
time, interpolation of the Sobolev scale gives
\begin{equation}
\label{eq:interpolacion-convergencia-C-Hsmenos1}
 \|\mathbf{w}\|_{H^{s-1}(M,\mathbf{E})}
 \leq C\|\mathbf{w}\|_{L^2(M,\mathbf{E})}^{\frac{1}{s}}
          \|\mathbf{w}\|_{H^s(M,\mathbf{E})}^{\frac{s-1}{s}},
\end{equation}
and, after integration in time,
\begin{equation}
\label{eq:interpolacion-convergencia-L2-Hs}
 \|\mathbf{w}\|_{L^2(0,T;H^s(M,\mathbf{E}))}
 \leq C\|\mathbf{w}\|_{L^2(0,T;H^1(M,\mathbf{E}))}^{\frac{1}{s}}
          \|\mathbf{w}\|_{L^2(0,T;H^{s+1}(M,\mathbf{E}))}^{\frac{s-1}{s}}.
\end{equation}
Thus $(\mathbf{u}^{(k)})$ converges strongly in
$C([0,T];H^{s-1}(M,\mathbf{E}))\cap L^2(0,T;H^s(M,\mathbf{E}))$. Weak compactness in $L^2(0,T;H^{s+1}(M,\mathbf E))$ and
$L^2(0,T;H^{s-1}(M,\mathbf E))$, applied to the sections and their time
derivatives, together with weak-$*$ compactness in
$L^\infty(0,T;H^s(M,\mathbf E))$, gives a subsequence with weak limits in
these spaces. The limit of the sections agrees with the strong limit
$\mathbf u$ already obtained. If $\mathbf z$ is the weak limit of the time
derivatives, integration by parts against smooth test sections compactly
supported in $(0,T)\times M$ passes to the limit and gives
$\mathbf z=\partial_t\mathbf u$. Thus $\mathbf u\in\mathcal X_T^s(\mathbf E)$.

Since $s-1>\displaystyle\frac{m}{2}+1$, convergence in $C([0,T];H^{s-1}(M,\mathbf E))$
implies uniform convergence of the first jets. These remain in a fixed
compact subset of $J^1(\mathcal U)$: their zeroth-order components lie in
$T_{\displaystyle\frac{\varepsilon}{2}}(\mathbf u_0)$ and their first derivatives are uniformly
bounded. Consequently,
\[
 A(\mathbf u^{(k)})\longrightarrow A(\mathbf u),
 \qquad
 F(j^1\mathbf u^{(k)})\longrightarrow F(j^1\mathbf u)
 \quad\text{uniformly on }M\times[0,T].
\]
For the principal term, writing $\|\cdot\|_{L^2_{t,x}}$ for the norm of
$L^2(0,T;L^2(M,T^*M\otimes\mathbf E))$, we have
\[
 \begin{aligned}
 &\|A(\mathbf u^{(k)})\nabla\mathbf u^{(k+1)}
                  -A(\mathbf u)\nabla\mathbf u\|_{L^2_{t,x}}\\
 &\quad\leq
 \|A(\mathbf u^{(k)})-A(\mathbf u)\|_\infty
                         \|\nabla\mathbf u^{(k+1)}\|_{L^2_{t,x}}\\
 &\qquad+\|A(\mathbf u)\|_\infty
                         \|\nabla(\mathbf u^{(k+1)}-\mathbf u)\|_{L^2_{t,x}}
 \longrightarrow0.
 \end{aligned}
\]
This uses convergence in $\mathcal Z_T(\mathbf E)$ and its uniform bound.
Equation \eqref{eq:iteracion-parabolica-cuaslineal}, tested against a smooth
compactly supported section, therefore passes to the limit and gives
\eqref{eq:ecuacion-parabolica-cuaslineal} distributionally. Both sides belong
to $L^2(0,T;H^{s-1}(M,\mathbf E))$, so the equality also holds in that space.
Convergence in $C([0,T];H^{s-1}(M,\mathbf E))$ preserves
$\mathbf u(0)=\mathbf u_0$ and the inclusion of the values in
$T_{\displaystyle\frac{\varepsilon}{2}}(\mathbf u_0)\subseteq\mathcal U$.

If $\mathbf{u}$ and $\mathbf{v}$ are two solutions, repeat the calculation
\eqref{eq:diferencia-A-iteracion-L2}--
\eqref{eq:energia-contraccion-L2-cuaslineal} for $\mathbf{w}=\mathbf{u}-\mathbf{v}$. Now
$\mathbf{w}(0)=0$, and the right-hand side contains $C_R\|\mathbf{w}(t)\|_{L^2(M,\mathbf{E})}^2$ in place of the
difference from the preceding step. Grönwall's lemma~\ref{lema: gronwall} gives $\mathbf{w}=0$ on all of $[0,T]$.
Uniformity of the existence time follows because all constants
used depend only on the bound $H^s(M,\mathbf E)$, distance from the complement
of $\mathcal U$, and the lower parabolicity bound.
\end{proof}

The autonomous form of the equation suffices for most of the chapter. Geometric
applications will also involve coefficients with explicit time
dependence; the same argument gives the following version.

\begin{corollary}[Explicit time dependence]
\label{cor:existencia-parabolica-cuaslineal-tiempo}
Let $T_*>0$ and let
\begin{align*}
 A&\colon[0,T_*]\times\mathcal U\longrightarrow
 \operatorname{End}(T^*M\otimes\mathbf{E}),\\
 F&\colon[0,T_*]\times J^1(\mathcal U)\longrightarrow\mathbf{E}
\end{align*}
be smooth bundle maps. Suppose that $A(t,z)$ is self-adjoint and the
strong parabolicity condition holds on $(0,\mathbf{u}_0)$ with a
positive constant. Then there exist $T\in(0,T_*]$ and a unique solution
$\mathbf{u}\in\mathcal X_T^s(\mathbf{E})$ of
\[
 \partial_t\mathbf{u}+\nabla^*\bigl(A(t,\mathbf{u})\nabla\mathbf{u}\bigr)
 =F(t,j^1\mathbf{u}),
 \qquad \mathbf{u}(0)=\mathbf{u}_0.
\]
The assertions on regularity, dependence on data, and uniformity of the
time remain valid when the bounds for derivatives of $A$ and $F$ are
uniform for $t\in[0,T_*]$ on the compact sets considered.
\end{corollary}

\begin{proof}
In the iteration of
\eqref{eq:iteracion-parabolica-cuaslineal}, replace the coefficients by
$A(t,\mathbf{u}^{(k)})$ and $F(t,j^1\mathbf{u}^{(k)})$. The composition estimates
of Lemma~\ref{lem:moser-parabolico-haces} are uniform in $t$
because the time interval is compact. The energy estimate,
remaining inside $\mathcal V$, and contraction in
$\mathcal Z_T(\mathbf{E})$ are unchanged. The proofs of regularity and parameter
dependence add only derivatives of the coefficients with respect to $t$, which
are uniformly bounded lower-order terms.
\end{proof}

\begin{corollary}[Regularity and dependence on data]
\label{cor:regularidad-dependencia-parabolica}
Let $(M,\mathbf{g}_0)$ be a closed Riemannian manifold, let
$\mathbf{E}\to M$ be a vector bundle with a bundle metric and compatible
connection, and consider the quasilinear equation
\eqref{eq:ecuacion-parabolica-cuaslineal}.
Under the hypotheses of
Theorem~\ref{teo:existencia-parabolica-cuaslineal}, the following assertions
hold.
\begin{enumerate}
 \item On bounded subsets of initial data sharing bounds for
 parabolicity and distance from $\mathbf{E}\setminus\mathcal U$, the solutions
 satisfy
 \begin{equation}
 \label{eq:dependencia-datos-cuaslineal}
 \|\mathbf{u}-\mathbf{v}\|_{\mathcal X_T^{s-1}(\mathbf{E})}
 \leq C\|\mathbf{u}_0-\mathbf{v}_0\|_{H^{s-1}(M,\mathbf{E})}.
 \end{equation}
 \item If $\mathbf{u}_0\in H^{s+q}(M,\mathbf{E})$ for some integer $q\geq0$, then, after
 reducing $T$ using only a bound for $\|\mathbf{u}_0\|_{H^s(M,\mathbf{E})}$, we have
 $\mathbf{u}\in\mathcal X_T^{s+q}(\mathbf{E})$.
 \item If $\mathbf{u}_0$ is smooth, the solution is smooth on $M\times[0,T]$. For
 merely $H^s(M,\mathbf E)$ initial data, the solution is smooth on
 $M\times(0,T]$.
 \item A smooth family of smooth initial data gives a smooth family
 of solutions, as long as they all remain in a set on which
 the preceding bounds are uniform.
\end{enumerate}
\end{corollary}

\begin{proof}
The equation for $\mathbf{u}-\mathbf{v}$ is
\begin{align*}
 \partial_t(\mathbf{u}-\mathbf{v})+\nabla^*(A(\mathbf{u})\nabla(\mathbf{u}-\mathbf{v}))
 ={}&-\nabla^*\bigl((A(\mathbf{u})-A(\mathbf{v}))\nabla \mathbf{v}\bigr)\\
 &+F(j^1\mathbf{u})-F(j^1\mathbf{v}).
\end{align*}
Apply $\nabla^q$, $0\leq q\leq s-1$, and repeat the integration by
parts in Lemma~\ref{lem:energia-parabolica-orden-alto}. For the term
containing the coefficient difference, use
\begin{equation}
\label{eq:producto-diferencia-soluciones-Hsmenos1}
 \|(A(\mathbf{u})-A(\mathbf{v}))\nabla \mathbf{v}\|_{H^{s-1}(M,\boldsymbol{\mathcal H})}
 \leq C_R
 \|\mathbf{u}-\mathbf{v}\|_{H^{s-1}(M,\mathbf{E})}\|\mathbf{v}\|_{H^s(M,\mathbf{E})}.
\end{equation}
This estimate is valid because $H^{s-1}(M)$ is an algebra. The difference of the
$F$ terms is bounded in $H^{s-2}(M,\mathbf{E})$ by
$C_R\|\mathbf{u}-\mathbf{v}\|_{H^{s-1}(M,\mathbf{E})}$. After absorbing the derivative of order $s$ of
$\mathbf{u}-\mathbf{v}$, we obtain
\begin{equation}
\label{eq:energia-dependencia-Hsmenos1}
 \frac{d}{dt}\mathscr E_{s-1}(\mathbf{u}-\mathbf{v})
 +c\|\mathbf{u}-\mathbf{v}\|_{H^s(M,\mathbf{E})}^2
 \leq C_R\mathscr E_{s-1}(\mathbf{u}-\mathbf{v}).
\end{equation}
Grönwall's lemma~\ref{lema: gronwall} controls the
$C([0,T];H^{s-1}(M,\mathbf{E}))\cap L^2(0,T;H^s(M,\mathbf{E}))$ part. The difference equation also controls
$\partial_t(\mathbf{u}-\mathbf{v})$ in $L^2(0,T;H^{s-2}(M,\mathbf{E}))$, giving
\eqref{eq:dependencia-datos-cuaslineal}.

We now prove persistence of regularity using an energy of arbitrary
order. Let $N\geq s$ and first suppose that $\mathbf{u}$ is smooth. On an
interval $[0,T]$ where $\displaystyle \|\mathbf{u}\|_{L^\infty(0,T;H^s(M,\mathbf{E}))}\leq R$, the Moser estimates with higher-order control
give
\begin{align}
 \|[\nabla^N,A(\mathbf{u})](\nabla \mathbf{u})\|_{L^2(M,T^{(0,N)}(TM)\otimes\boldsymbol{\mathcal H})}
 +\|F(j^1\mathbf{u})\|_{H^{N-1}(M,\mathbf{E})}
 \leq C_{N,R}\bigl(1+\|\mathbf{u}\|_{H^N(M,\mathbf{E})}\bigr).
 \label{eq:moser-alto-orden-cuaslineal}
\end{align}
Indeed, \eqref{eq:moser-conmutador-parabolico} gives
\[
\begin{aligned}
 &\|[\nabla^N,A(\mathbf{u})]\nabla\mathbf{u}\|_
       {L^2(M,T^{(0,N)}(TM)\otimes\boldsymbol{\mathcal H})}\\
 &\quad\leq C_N
 \|\nabla A(\mathbf{u})\|_
       {L^\infty(M,T^*M\otimes\operatorname{End}(\boldsymbol{\mathcal H}))}
 \|\nabla\mathbf{u}\|_{H^{N-1}(M,\boldsymbol{\mathcal H})}\\
 &\qquad+C_N
 \|A(\mathbf{u})\|_{H^N(M,\operatorname{End}(\boldsymbol{\mathcal H}))}
 \|\nabla\mathbf{u}\|_{L^\infty(M,\boldsymbol{\mathcal H})}.
\end{aligned}
\]
The bound $H^s(M,\mathbf E)\hookrightarrow C^2(M,\mathbf E)$ controls both factors with
infinity norms. Estimate \eqref{eq:moser-composicion-lineal-orden-alto}
controls $A(\mathbf{u})$ in $H^N(M,\operatorname{End}(T^*M\otimes\mathbf E))$ by
$C_{N,R}(1+\|\mathbf{u}\|_{H^N(M,\mathbf E)})$. Applied to
$(\mathbf{u},\nabla\mathbf{u})\mapsto F(\mathbf{u},\nabla\mathbf{u})$ at
order $N-1$, the same estimate controls the remaining term in
\eqref{eq:moser-alto-orden-cuaslineal}. All constants use the bounds
for the first jets and the fixed structures.
Applying $\nabla^N$ to the
equation, integrating by parts, and using
\eqref{eq:moser-alto-orden-cuaslineal} yields
\begin{equation}
\label{eq:energia-alto-orden-cuaslineal}
 \frac{d}{dt}\mathscr E_N(\mathbf{u}(t))
 +c\|\mathbf{u}(t)\|_{H^{N+1}(M,\mathbf{E})}^2
 \leq C_{N,R}\bigl(1+\mathscr E_N(\mathbf{u}(t))\bigr).
\end{equation}
For a solution $\mathbf u\in\mathcal X_T^N(\mathbf E)$,
the identity allowing this energy to be used is
\begin{equation}
\label{eq:energia-covariante-tiempo-debil}
 \frac{d}{dt}\mathscr E_N(\mathbf u(t))
 =2\mathcal B_N(\partial_t\mathbf u(t),\mathbf u(t))
 \quad\text{for almost every }t\in(0,T),
\end{equation}
where $\mathcal B_N$ is the form extended in the proof of
Lemma~\ref{lem:energia-parabolica-orden-alto}.
To verify it, fix a compact interval $J\subset(0,T)$ and
regularize only in time:
$\mathbf u_\varepsilon=\rho_\varepsilon*_t\mathbf u$, with
$\varepsilon$ smaller than the distance from $J$ to the endpoints.
Continuity of translations in Bochner spaces gives
\[
 \mathbf u_\varepsilon\to\mathbf u
 \quad\text{in }L^2(J;H^{N+1}(M,\mathbf E)),\qquad
 \partial_t\mathbf u_\varepsilon\to\partial_t\mathbf u
 \quad\text{in }L^2(J;H^{N-1}(M,\mathbf E)).
\]
The first convergence implies
$\mathscr E_N(\mathbf u_\varepsilon)\to\mathscr E_N(\mathbf u)$
in $L^1(J)$. By the bound for $\mathcal B_N$, the error on the right-hand
side, integrated over $J$, is at most
\begin{align*}
 &C\|\partial_t\mathbf u_\varepsilon-\partial_t\mathbf u\|_
       {L^2(J;H^{N-1}(M,\mathbf E))}
    \|\mathbf u_\varepsilon\|_{L^2(J;H^{N+1}(M,\mathbf E))}\\
 &\quad+C\|\partial_t\mathbf u\|_{L^2(J;H^{N-1}(M,\mathbf E))}
    \|\mathbf u_\varepsilon-\mathbf u\|_{L^2(J;H^{N+1}(M,\mathbf E))},
\end{align*}
which tends to zero. For $\mathbf u_\varepsilon$, differentiation of
the covariant norm and the formula for $\mathcal B_N$ give the classical
identity; testing against a function in $C_c^\infty(J)$ and taking
limits gives \eqref{eq:energia-covariante-tiempo-debil} in
distributions. Its right-hand side belongs to $L^1(0,T)$.
Lemma~\ref{lem:traza-temporal-parabolica} provides the continuous
representative in $H^N(M,\mathbf E)$, so the energy has the
absolutely continuous representative determined by this identity.

Finally, for almost every $t$, substitute the equation into
$\mathcal B_N(\partial_t\mathbf u,\mathbf u)$.
Spatial integrations by parts are valid for
$\mathbf u(t)\in H^{N+1}(M,\mathbf E)$ by the explicit definition of
that form. The product and Moser estimates in
\eqref{eq:moser-alto-orden-cuaslineal} control the resulting
terms in the same spaces; they are also obtained as limits
of smooth approximations in $H^{N+1}(M,\mathbf E)$.
Thus the energy estimate holds for almost every time and can be
integrated on any subinterval of $[0,T]$.

Let us justify application of this estimate on the existence interval
determined by the norm of order $s$. If $\mathbf{u}_0\in H^N(M,\mathbf E)$, the quasilinear
theorem applied initially at order $N$ gives a solution on a
possibly smaller interval. By uniqueness, it agrees there with the solution
of order $s$. While both are defined, Grönwall applied to
\eqref{eq:energia-alto-orden-cuaslineal} bounds its $H^N(M,\mathbf E)$ norm in terms of
$N,R,T$ and $\|\mathbf{u}_0\|_{H^N(M,\mathbf E)}$. Distance to the boundary of
$\mathcal U$ and parabolicity retain the bounds from the order-$s$
construction. Restarting the existence theorem at any time in
this interval, these bounds and the bound $H^N(M,\mathbf E)$ give a positive
additional time independent of the restart time. If the interval of
$H^N(M,\mathbf E)$ regularity ended before $T$, restarting sufficiently
near its endpoint would extend beyond it. By uniqueness, the extension would
remain the same solution. This gives $H^N(M,\mathbf E)$ regularity up to $T$.

To justify the identities for $H^N(M,\mathbf E)$ data, approximate the data in
that norm by smooth sections. The order-$s$ estimate gives
a common interval, and the preceding argument and
\eqref{eq:energia-alto-orden-cuaslineal} give uniform bounds in
$\mathcal X_T^N$. Weak passage to the limit and identification of the trace
are as in the linear theorem; uniqueness identifies the limit with
$\mathbf{u}$. We conclude that $\mathbf{u}\in\mathcal X_T^N$. The interval $T$
is fixed using the order-$s$ bounds, although the order-$N$ norm of the
solution also depends on that of the initial data.

The same energy proves smoothing without attributing to the data
derivatives they do not have. Approximate \(\mathbf u_0\) in \(H^s\)
by smooth data \(\mathbf u_{0,l}\), and let \(\mathbf u_l\) be their
solutions. Persistence of regularity already proved and the order-\(s\)
bounds give a common interval \([0,T]\), smooth solutions, and
\[
 \sup_l\|\mathbf u_l\|_{\mathcal X_T^s}^2
 \leq C_R(1+\|\mathbf u_0\|_{H^s}^2).
\]
All images remain in the same uniformly parabolic compact set.
Fix \(0<\tau<T\), \(q\geq1\), and set
\(\tau_j:=\displaystyle\frac{j\tau}{q}\) for \(0\leq j\leq q\).
Choose cutoffs \(0\leq\chi_j\leq1\) vanishing near
\([0,\tau_{j-1}]\), equal to one on \([\tau_j,T]\), and satisfying
\(\|\chi_j'\|_\infty\leq C_\chi\displaystyle\frac{q}{\tau}\),
obtained by translation and dilation of a single cutoff.

For \(N=s+j\), write
\(Y_{j,l}(t):=\chi_j(t)^2\mathscr E_N(\mathbf u_l(t))\).
Multiplying \eqref{eq:energia-alto-orden-cuaslineal} by
\(\chi_j^2\) and differentiating the product gives
\[
 \frac{d}{dt}Y_{j,l}
 +c\chi_j^2\|\mathbf u_l\|_{H^{N+1}}^2
 \leq C_{N,R}(1+Y_{j,l})
       +2\|\chi_j\chi_j'\|_\infty
           \mathbf 1_{(\tau_{j-1},\tau_j)}\mathscr E_N(\mathbf u_l).
\]
Also \(Y_{j,l}(0)=0\). Grönwall and equivalence of energy with the
Sobolev norm give, with constants independent of \(l\),
\[
 \begin{aligned}
 &\sup_{\tau_j\leq t\leq T}\|\mathbf u_l(t)\|_{H^N}^2
            +\int_{\tau_j}^T\|\mathbf u_l(t)\|_{H^{N+1}}^2\,dt\\
 &\quad\leq C_{j,q,R,\tau}
       \left(1+\int_{\tau_{j-1}}^{\tau_j}
                                  \|\mathbf u_l(t)\|_{H^N}^2\,dt\right).
 \end{aligned}
\]
For \(j=1\), the common \(\mathcal X_T^s\) bounds control the
integral on the right. If the bound has been obtained for \(j-1\),
its integral term controls
\(\int_{\tau_{j-1}}^T\|\mathbf u_l\|_{H^{s+j}}^2\) and hence the
right-hand side at step \(j\). This induction gives uniform bounds
for every \(j\leq q\), without using high norms of the approximating
data. The equation and Moser estimates also give
\[
 \|\partial_t\mathbf u_l\|_{L^2(\tau_j,T;H^{N-1})}
 \leq C_{N,R}\left(1+
              \|\mathbf u_l\|_{L^2(\tau_j,T;H^{N+1})}\right).
\]
Weak passage to the limit, identification of the equation, and the
trace lemma thus give order-\(s+j\) regularity of the limit on
\([\tau_j,T]\). Semicontinuity preserves the bound
\begin{align}
 &\|\mathbf u\|_{L^\infty(\tau_j,T;H^{s+j}(M,\mathbf E))}^2
       +\int_{\tau_j}^T\|\mathbf u(t)\|_{H^{s+j+1}(M,\mathbf E)}^2\,dt
 \nonumber\\
 &\qquad\leq C_{j,q,R,\tau}(1+\|\mathbf u_0\|_{H^s(M,\mathbf E)}^2).
 \label{eq:suavizado-parabolico-tiempo-positivo}
\end{align}

The induction reaches every $q$. Sobolev embedding gives spatial
smoothness on $M\times(0,T]$, and the equation expresses every time derivative in
terms of spatial derivatives already controlled; this gives
joint smoothness. If $\mathbf{u}_0$ is smooth, estimate
\eqref{eq:energia-alto-orden-cuaslineal} starts at $t=0$ for every $N$, and the
solution is smooth up to the initial time.

Finally, we justify smooth dependence. Let
\(a\mapsto\mathbf u_{0,a}\) be a smooth family of smooth data.
Work on a compact parameter interval and a common time interval.
Persistence of regularity gives uniform bounds for the solutions
\(\mathbf u_a\) in every \(\mathcal X_T^N\).
The difference energy already proved, applied at order \(r\) with
uniform order-\(r+1\) bounds, gives, for every \(r\geq s\),
\[
 \|\mathbf u_b-\mathbf u_a\|_{\mathcal X_T^r}
 \leq C_r\|\mathbf u_{0,b}-\mathbf u_{0,a}\|_{H^r}.
\]
Constants are taken locally in the parameter. In particular, solutions
depend continuously on it at all orders.

Fix \(a\), write \(u:=\mathbf u_a\), \(v:=\mathbf u_{a+\varepsilon}\),
and set
\[
 h_\varepsilon:=\frac{v-u}{\varepsilon},\qquad
 u_\theta:=u+\theta(v-u),\qquad 0\leq\theta\leq1.
\]
Use the smooth extensions of \(A,F\) to the entire fibers constructed
in the theorem's proof. Define the linear coefficients
\[
 \begin{aligned}
 \mathcal A_\varepsilon[h]&:=\int_0^1D_zA(u_\theta)[h]\,d\theta,\\
 \mathcal B_\varepsilon[h]&:=\int_0^1D_zF(j^1u_\theta)[h]\,d\theta,\\
 \mathcal C_\varepsilon[p]&:=\int_0^1D_pF(j^1u_\theta)[p]\,d\theta,\\
 \mathcal D_\varepsilon[h]&:=\mathcal A_\varepsilon[h]\nabla v.
 \end{aligned}
\]
The quotient satisfies exactly
\[
 \partial_t h_\varepsilon+
 \nabla^*(A(u)\nabla h_\varepsilon+\mathcal D_\varepsilon[h_\varepsilon])
 =\mathcal B_\varepsilon[h_\varepsilon]
       +\mathcal C_\varepsilon[\nabla h_\varepsilon],
 \qquad h_\varepsilon(0)=
             \frac{\mathbf u_{0,a+\varepsilon}-\mathbf u_{0,a}}{\varepsilon}.
\]
Expanding \(\nabla^*(\mathcal D_\varepsilon[h])\) produces first-order
and zeroth-order terms, leaving the principal part unchanged.
Denote the limiting coefficients by
\[
 \mathcal D[h]:=D_zA(u)[h]\nabla u,\qquad
 \mathcal B:=D_zF(j^1u),\qquad \mathcal C:=D_pF(j^1u).
\]
Uniform higher-order bounds and Moser estimates imply, for \(r\geq s\),
\[
 \begin{aligned}
 d_{r,\varepsilon}:={}&
 \|\mathcal D_\varepsilon-\mathcal D\|_{C([0,T];H^{r+1})}
 +\|\mathcal B_\varepsilon-\mathcal B\|_{C([0,T];H^r)}
 +\|\mathcal C_\varepsilon-\mathcal C\|_{C([0,T];H^r)}\\
 \leq{}&C_r\|v-u\|_{C([0,T];H^{r+2})}
 \leq C_r|\varepsilon|
             \sup_b\|\partial_a\mathbf u_{0,b}\|_{H^{r+2}}.
 \end{aligned}
\]
Coefficient norms are taken in their homomorphism bundles;
parameter suprema range over the fixed compact set. The last inequality
uses the difference estimate at order \(r+2\) and the fundamental
theorem of calculus for the data.

Linear theory constructs a unique solution \(h\), with data
\(h(0)=\partial_a\mathbf u_{0,a}\), of
\begin{align}
 \partial_t h+\nabla^*\bigl(&A(u)\nabla h+D_zA(u)[h]\nabla u\bigr)
 \nonumber\\
 &={}D_zF(u,\nabla u)[h]+D_pF(u,\nabla u)[\nabla h].
 \label{eq:ecuacion-variacion-parametrica-parabolica}
\end{align}
Its regularity holds at every order. Subtracting its equation from the
quotient equation and applying the linear estimate with the quotient's
coefficients, product estimates give
\[
 \begin{aligned}
 \|h_\varepsilon-h\|_{\mathcal X_T^r}
 \leq C_r\biggl(&
 \left\|\frac{\mathbf u_{0,a+\varepsilon}-\mathbf u_{0,a}}{\varepsilon}
                            -\partial_a\mathbf u_{0,a}\right\|_{H^r}
       +d_{r,\varepsilon}\|h\|_{\mathcal X_T^r}\biggr).
 \end{aligned}
\]
Indeed, the difference of \(\mathcal D\) contributes
\(\nabla^*((\mathcal D_\varepsilon-\mathcal D)[h])\) in
\(L^2(0,T;H^{r-1})\); differences of \(\mathcal B,\mathcal C\)
belong to the same space, with norms at most
\(C_rd_{r,\varepsilon}\|h\|_{\mathcal X_T^r}\).
For both signs of \(\varepsilon\), the integral remainder of the data gives
\[
 \left\|\frac{\mathbf u_{0,a+\varepsilon}-\mathbf u_{0,a}}{\varepsilon}
                            -\partial_a\mathbf u_{0,a}\right\|_{H^r}
 \leq\frac{|\varepsilon|}{2}
                 \sup_b\|\partial_a^2\mathbf u_{0,b}\|_{H^r}.
\]
Thus \(\|h_\varepsilon-h\|_{\mathcal X_T^r}\to0\), with constants
independent of \(\varepsilon\); this identifies
\(h=\partial_a\mathbf u_a\) at every order. Comparing the linearized
equations at nearby parameters controls the difference of their
solutions by differences of initial data and coefficients in the
same norms. Both terms tend to zero, so the derivative constructed
is continuous in the parameter.

Let us make the passage to higher parameter derivatives explicit. For
a real variable $a$ write $\mathbf{u}^{[k]}=\partial_a^k\mathbf{u}_a$
and define the nonlinear differential operator
\[
 \mathcal H(\mathbf{v}):=-\nabla^*(A(\mathbf{v})\nabla\mathbf{v})
                          +F(\mathbf{v},\nabla\mathbf{v}).
\]
Suppose that the derivatives
$\mathbf{u}^{[1]},\ldots,\mathbf{u}^{[q-1]}$ have been constructed and are continuous at all spatial orders,
with $q\geq2$. Formula \eqref{eq:faa-di-bruno-curva-coeficientes}, in the version
of Corollary~\ref{cor:faa-di-bruno-espacios-normados}, gives the
candidate equation
\[
 (\partial_t-D\mathcal H(\mathbf{u}_a))\mathbf{u}^{[q]}
 =\sum_{j=2}^q\ \sum_{\substack{k_1,\ldots,k_j\in\mathbb N\\
                                  k_1+\cdots+k_j=q}}
 \frac{q!}{j!\,k_1!\cdots k_j!}
 D^j\mathcal H(\mathbf{u}_a)
 [\mathbf{u}^{[k_1]},\ldots,\mathbf{u}^{[k_j]}].
\]
In this sum, every $k_i$ is less than $q$; by the induction hypothesis, the
right-hand side is known. The operator on the left is precisely
that of \eqref{eq:ecuacion-variacion-parametrica-parabolica} and retains the
principal part $\nabla^*(A(\mathbf{u}_a)\nabla\cdot)$. The initial data
are $\partial_a^q\mathbf{u}_{0,a}$. Linear theory gives a unique
solution at every spatial order.

To identify it as a derivative, let
\[
 h_{q,\varepsilon}:=
  \frac{\mathbf u^{[q-1]}_{a+\varepsilon}-\mathbf u^{[q-1]}_a}{\varepsilon}.
\]
Denote by \(R_{q-1,a}\) the known source in the equation for
\(\mathbf u_a^{[q-1]}\), with \(R_{1,a}:=0\). Subtracting the equations
at \(a+\varepsilon\) and \(a\) gives exactly
\[
 (\partial_t-D\mathcal H(v))h_{q,\varepsilon}=R_{q,\varepsilon},
\]
where the integral composition rule yields
\[
 \begin{aligned}
 R_{q,\varepsilon}:={}&
 \int_0^1D^2\mathcal H(u_\theta)
              [h_\varepsilon,\mathbf u_a^{[q-1]}]\,d\theta
 +\frac{R_{q-1,a+\varepsilon}-R_{q-1,a}}{\varepsilon}.
 \end{aligned}
\]
The source \(R_{q-1,a}\) depends only on \(\mathbf u_a\) and its
parameter derivatives up to order \(q-2\). By the induction hypothesis
it is continuously differentiable in \(\mathcal Y_T^r\), using two
additional spatial orders to control the operator's derivatives.
Its quotient is therefore
\(\int_0^1\partial_aR_{q-1,a+\theta\varepsilon}\,d\theta\).
The limiting source is
\[
 R_q=D^2\mathcal H(u)
       [\mathbf u_a^{[1]},\mathbf u_a^{[q-1]}]
                  +\partial_aR_{q-1,a},
\]
which agrees with the right-hand side of the candidate equation by
the product rule and Faà di Bruno. To quantify convergence, set
\(\mathbf u_b^{[0]}:=\mathbf u_b\) and
\[
 \omega_{k,r+2}(\varepsilon):=
 \sup_{|b-a|\leq|\varepsilon|}
 \|\mathbf u_b^{[k]}-\mathbf u_a^{[k]}\|_{C([0,T];H^{r+2})},
 \qquad 0\leq k\leq q-1.
\]
Moser estimates, the preceding integral formula, and the uniform
bounds on the derivatives already constructed give
\[
 \|R_{q,\varepsilon}-R_q\|_{\mathcal Y_T^r}
 \leq C_{q,r}\left(
 \|h_\varepsilon-\mathbf u_a^{[1]}\|_{\mathcal X_T^{r+2}}
            +\sum_{k=0}^{q-1}\omega_{k,r+2}(\varepsilon)\right).
\]
The first term tends to zero by the first-derivative case,
and each modulus \(\omega_{k,r+2}\) tends to zero by the continuity
already proved. The constant is independent of \(\varepsilon\).
Let \(d^{(q)}_{r,\varepsilon}\) be the sum of the differences between
the coefficients of \(D\mathcal H(v)\) and \(D\mathcal H(u)\), in
\(C([0,T];H^r)\), with one additional spatial derivative for those
that remain inside a divergence. The same estimates give
\[
 d^{(q)}_{r,\varepsilon}
 \leq C_r\|v-u\|_{C([0,T];H^{r+2})}
 \leq C_r|\varepsilon|\sup_b
                         \|\partial_a\mathbf u_{0,b}\|_{H^{r+2}},
\]
so this difference also tends to zero.
Subtracting the candidate equation and applying the same linear estimate gives
\[
 \begin{aligned}
 \|h_{q,\varepsilon}-\mathbf u_a^{[q]}\|_{\mathcal X_T^r}
 \leq C_r\biggl(&
 \left\|\frac{\partial_a^{q-1}\mathbf u_{0,a+\varepsilon}
                 -\partial_a^{q-1}\mathbf u_{0,a}}{\varepsilon}
                     -\partial_a^q\mathbf u_{0,a}\right\|_{H^r}\\
 &+\|R_{q,\varepsilon}-R_q\|_{\mathcal Y_T^r}
   +d^{(q)}_{r,\varepsilon}\|\mathbf u_a^{[q]}\|_{\mathcal X_T^r}\biggr).
 \end{aligned}
\]
The first term is at most
\(\displaystyle\frac{|\varepsilon|}{2}
       \sup_b\|\partial_a^{q+1}\mathbf u_{0,b}\|_{H^r}\),
and the other two tend to zero by the stated convergences.
Constants are uniform for small \(\varepsilon\) on the fixed compact
parameter set. Thus the quotient converges to the candidate linear
solution. Comparing two parameter values with the same estimate
proves continuity of this derivative and closes the induction.
For parameters in an open subset of \(\mathbb R^d\), apply the
argument to successive directional quotients; the multilinear
Faà di Bruno formulas give the mixed derivatives.

This proves the last
assertion.
\end{proof}

\begin{theorem}[Continuation criterion]
\label{teo:criterio-continuacion-parabolico}
Let $(M,\mathbf{g}_0)$ be a closed Riemannian manifold, let
$\mathbf{E}\to M$ be a vector bundle with a bundle metric and compatible
connection, and let $\mathbf{u}$ be a maximal solution of
\eqref{eq:ecuacion-parabolica-cuaslineal} defined on
$[0,T_{\displaystyle\max})$, with $T_{\displaystyle\max}<\infty$, and suppose that
$\mathbf{u}\in\mathcal X_T^s(\mathbf{E})$ for every $T<T_{\displaystyle\max}$. If there exist $R,\delta>0$
such that,
for every $0\leq t<T_{\displaystyle\max}$,
\begin{align}
 \|\mathbf{u}(t)\|_{H^s(M,\mathbf{E})}&\leq R,
 \label{eq:continuacion-cota-sobolev}\\
 \inf_{x\in M}\operatorname{dist}_{\mathbf{E}_x}
 \bigl(\mathbf{u}(x,t),\mathbf{E}_x\setminus\mathcal U_x\bigr)&\geq\delta,
 \label{eq:continuacion-distancia-abierto}\\
 \left\langle A(\mathbf{u}(x,t))\zeta,\zeta\right\rangle_{\mathbf{g}_0,\mathbf{h}_{\mathbf{E}}}
 &\geq\delta|\zeta|_{\mathbf{g}_0,\mathbf{h}_{\mathbf{E}}}^2.
 \label{eq:continuacion-parabolicidad}
\end{align}
The last inequality is required for every
$x\in M$ and $\zeta\in T_x^*M\otimes \mathbf{E}_x$. Then $\mathbf{u}$ extends to an
interval $[0,T_{\displaystyle\max}+\varepsilon)$, contradicting
maximality. Thus, if the maximal time is finite,
at least one quantity controlled in
\eqref{eq:continuacion-cota-sobolev}--\eqref{eq:continuacion-parabolicidad}
must degenerate.
\end{theorem}

\begin{proof}
The last assertion of
Theorem~\ref{teo:existencia-parabolica-cuaslineal} gives a time
$\varepsilon>0$ depending only on $R$, $\delta$, and bounds for $A,F$ on
the compact set traversed by the first jets of $\mathbf{u}$. These latter bounds
follow from $H^s(M,\mathbf{E})\hookrightarrow C^1(M,\mathbf{E})$. For any
$t_0\in(T_{\displaystyle\max}-\displaystyle\frac{\varepsilon}{2},T_{\displaystyle\max})$, apply the existence
theorem with initial data $\mathbf{u}(t_0)$. This gives a solution on
$[t_0,t_0+\varepsilon]$. Uniqueness makes it agree with $\mathbf{u}$ on
$[t_0,T_{\displaystyle\max})$, thus extending it beyond $T_{\displaystyle\max}$.
\end{proof}

\section{Maximum principles}

Energy estimates control the solution in global norms, but do not
by themselves record the sign of a function or whether a tensor belongs to
a convex cone. Maximum principles provide this pointwise
information and will be the basic tool for studying curvature under
Ricci flow.

\begin{theorem}[Scalar maximum principle]
\label{teo:principio-maximo-parabolico-escalar}
Let $(M,\mathbf{g}_0)$ be a closed Riemannian manifold. Let
\[
 L_t\phi=\sum_{i,j=1}^m a^{ij}(x,t)(\nabla^2\phi)_{ij}
          +\sum_{i=1}^m b^i(x,t)\nabla_i\phi+c(x,t)\phi,
\]
where $a^{ij}$, $b^i$, and $c$ are continuous on $M\times[0,T]$,
$(a^{ij})$ is symmetric, and there exists $\lambda>0$ such that
$\displaystyle \sum_{i,j=1}^m a^{ij}(x,t)\xi_i\xi_j\geq\lambda|\xi|_{\mathbf{g}_0}^2$
for every $(x,t)\in M\times[0,T]$ and every $\xi\in T_x^*M$. Suppose that
$f\in C^0(M\times[0,T])$ and that
$u\in C^0(M\times[0,T])\cap C^{2,1}(M\times(0,T])$ satisfies
\begin{equation}
\label{eq:hipotesis-principio-maximo-escalar}
 \partial_tu-L_tu\leq f
\end{equation}
and let $\gamma\in\mathbb R$ be an upper bound for $c$ on
$M\times[0,T]$. Then, for $0\leq t\leq T$,
\begin{align}
 \sup_{x\in M}u(x,t)
 \leq e^{\gamma t}\biggl(&
 \max\left\{0,\sup_{x\in M}u(x,0)\right\}+\int_0^t e^{-\gamma\rho}
 \max\left\{0,\sup_{x\in M}f(x,\rho)\right\}\,d\rho
 \biggr).
 \label{eq:estimacion-principio-maximo-escalar}
\end{align}
In particular, if $c\leq0$, $f\leq0$, and $u(0)\leq0$, then $u\leq0$ on
$M\times[0,T]$.
\end{theorem}

\begin{proof}
Set $v(x,t)=e^{-\gamma t}u(x,t)$ and
\[
 F_+(t):=\max\left\{0,\sup_{x\in M}e^{-\gamma t}f(x,t)\right\},
 \qquad
 z(t):=\max\left\{0,\sup_{x\in M}u(x,0)\right\}
       +\int_0^tF_+(\rho)\,d\rho.
\]
The function $w=v-z$ satisfies
\begin{equation}
\label{eq:reduccion-principio-maximo-escalar}
 \partial_tw-\sum_{i,j=1}^m a^{ij}(\nabla^2w)_{ij}-\sum_{i=1}^m b^i\nabla_iw-(c-\gamma)w\leq0,
 \qquad w(0)\leq0,
\end{equation}
because $c-\gamma\leq0$ and $z\geq0$.

If $w$ attained a positive value, then for small $\varepsilon>0$ the function
$w_\varepsilon=w-\varepsilon t-\varepsilon$ would have a first nonnegative
maximum at some $(x_0,t_0)$ with $t_0>0$. There,
\[
 \partial_tw_\varepsilon\geq0,
 \qquad \nabla w_\varepsilon=0,
 \qquad \sum_{i,j=1}^m a^{ij}(\nabla^2w_\varepsilon)_{ij}\leq0.
\]
If $\displaystyle D=\partial_t-\displaystyle\sum_{i,j=1}^m a^{ij}(\nabla^2\,\cdot\,)_{ij}-\displaystyle\sum_{i=1}^m b^i\nabla_i-(c-\gamma)$, the three
preceding relations and $c-\gamma\leq0$ imply
$Dw_\varepsilon(x_0,t_0)\geq0$. On the other hand,
\[
 Dw_\varepsilon
 =Dw-\varepsilon+(c-\gamma)\varepsilon(t+1)
 \leq-\varepsilon,
\]
a contradiction. Therefore,
$w\leq0$. Multiplying by $e^{\gamma t}$ gives
\eqref{eq:estimacion-principio-maximo-escalar}.
\end{proof}

\begin{corollary}[Scalar comparison]
\label{cor:comparacion-parabolica-escalar}
Let $(M,\mathbf{g}_0)$ be a closed Riemannian manifold and
let $L_t$ be the operator of Theorem~\ref{teo:principio-maximo-parabolico-escalar}.
Let $u,v\in C^0(M\times[0,T])\cap C^{2,1}(M\times(0,T])$. If
\[
 \partial_tu-L_tu\leq\partial_tv-L_tv
 \quad\hbox{and}\quad u(0)\leq v(0),
\]
then $u\leq v$ on $M\times[0,T]$ when $c\leq0$. For bounded $c$, the
same conclusion follows after multiplying both functions by
$e^{-\gamma t}$ with $\displaystyle \gamma\geq\displaystyle\sup_{(x,t)\in M\times[0,T]}c(x,t)$.
\end{corollary}

\begin{proof}
Apply Theorem~\ref{teo:principio-maximo-parabolico-escalar} to $u-v$.
\end{proof}

To formulate the tensor version, a family
$\mathcal K=(K_x)_{x\in M}$ of subsets of $\mathbf{E}_x$ is called
\emph{invariant under parallel transport} if transport along
any curve takes $K_x$ into the set corresponding to its endpoint. In
particular, this condition identifies all the $K_x$ locally using the
connection. If $q\in\partial K_x$, a vector $\nu\in \mathbf{E}_x$ is an outward supporting
normal if
\[
 \langle\nu,z-q\rangle_{\mathbf{E}_x}\leq0
 \qquad\text{for every }z\in K_x.
\]

\begin{theorem}[Tensor maximum principle]
\label{teo:principio-maximo-parabolico-tensorial}
Let $(M,\mathbf{g}_0)$ be a closed Riemannian manifold and
let $\mathbf{E}\to M$ be a real vector bundle with a bundle metric and compatible
connection.
Let $\mathcal K=(K_x)_{x\in M}$ be a family of nonempty closed convex
sets invariant under parallel transport. Let
$\mathbf{S}\in C^0(M\times[0,T],\mathbf{E})\cap C^{2,1}(M\times(0,T],\mathbf{E})$ be a solution of
\begin{equation}
\label{eq:sistema-principio-maximo-tensorial}
 \partial_t\mathbf{S}=\sum_{i,j=1}^m a^{ij}(x,t)(\nabla^2\mathbf{S})_{ij}
              +\sum_{i=1}^m b^i(x,t)\nabla_i\mathbf{S}+\Phi(x,t,\mathbf{S}),
\end{equation}
where $\mathbf{a}\in C^0(M\times[0,T],S^2TM)$ and
$\mathbf{b}\in C^0(M\times[0,T],TM)$ are fields whose components
$a^{ij}$ and $b^i$ multiply the
identity of $\mathbf{E}$, and $(a^{ij})$ is symmetric and uniformly positive.
Suppose that $\Phi$ is continuous and locally
Lipschitz in its last variable, uniformly in
$(x,t)$ on compact sets, and that for every $q\in\partial K_x$ and every
outward supporting normal $\nu$ at $q$ we have
\begin{equation}
\label{eq:condicion-tangencia-convexa-parabolica}
 \langle\nu,\Phi(x,t,q)\rangle_{\mathbf{h}_{\mathbf{E}}}\leq0.
\end{equation}
If $\mathbf{S}(x,0)\in K_x$ for every $x\in M$, then
$\mathbf{S}(x,t)\in K_x$ for every $(x,t)\in M\times[0,T]$.
\end{theorem}

\begin{proof}
Let $\Pi_x\colon \mathbf{E}_x\longrightarrow K_x$ be the metric projection, unique
by convexity. Define
\[
 \rho(x,t):=\frac{1}{2}
 \left|\mathbf{S}(x,t)-\Pi_x(\mathbf{S}(x,t))\right|_{\mathbf{E}_x}^2.
\]
The squared-distance function to a closed convex set in a Euclidean
space is of class $C^1$ and has gradient
$z\mapsto z-\Pi_xz$. It is also convex and of class $C^{1,1}$. Invariance
of $\mathcal K$ under parallel transport identifies the distance
functions on neighboring fibers. In particular, $\rho$ is continuous, is of class
$C^1$ in positive time, and
\begin{equation}
\label{eq:derivada-temporal-distancia-convexo}
 \partial_t\rho(x,t)
 =\langle \mathbf{S}(x,t)-\Pi_x\mathbf{S}(x,t),\partial_t\mathbf{S}(x,t)\rangle_{\mathbf{h}_{\mathbf{E}}}.
\end{equation}

Fix $t>0$ and a point $x$ where
$\rho(\cdot,t)$ attains its maximum. Set
$q=\Pi_x\mathbf{S}(x,t)$ and $\nu=\mathbf{S}(x,t)-q$. Let $X\in T_xM$, let
$\gamma(h)=\exp_x(hX)$, and parallel-transport $\mathbf{S}(\gamma(h),t)$
back to $\mathbf{E}_x$; denote the transported curve by $z(h)$. Transport invariance
identifies $K_{\gamma(h)}$ with $K_x$, and maximality of $x$
gives
\[
 \frac{1}{2}\operatorname{dist}_{\mathbf{E}_x}(z(h),K_x)^2
 \leq\frac{1}{2}\operatorname{dist}_{\mathbf{E}_x}(z(0),K_x)^2.
\]
The subgradient inequality for the convex squared-distance function
is
\[
 \frac{1}{2}\operatorname{dist}_{\mathbf{E}_x}(z,K_x)^2
 \geq\frac{1}{2}\operatorname{dist}_{\mathbf{E}_x}(z(0),K_x)^2
      +\langle\nu,z-z(0)\rangle_{\mathbf{h}_{\mathbf{E}}}.
\]
Apply it to $z(h)$ and $z(-h)$, add, and divide by $h^2$. Letting
$h\to0$ gives
\begin{equation}
\label{eq:segunda-derivada-apoyo-convexo}
 \langle\nu,(\nabla^2\mathbf{S})(X,X)(x,t)\rangle_{\mathbf{h}_{\mathbf{E}}}\leq0.
\end{equation}
The first derivative at $h=0$ also gives
\begin{equation}
\label{eq:primera-derivada-apoyo-convexo}
 \langle\nu,\nabla_X\mathbf{S}(x,t)\rangle_{\mathbf{h}_{\mathbf{E}}}=0.
\end{equation}
Diagonalize $(a^{ij}(x,t))$ in an orthonormal basis of $T_xM$. Summing
\eqref{eq:segunda-derivada-apoyo-convexo} with the positive eigenvalues as
weights and using~\eqref{eq:primera-derivada-apoyo-convexo} yields
\begin{equation}
\label{eq:parte-espacial-apoyo-convexo}
 \left\langle\nu,
 \sum_{i,j=1}^m a^{ij}(\nabla^2\mathbf{S})_{ij}+\sum_{i=1}^m b^i\nabla_i\mathbf{S}\right\rangle_{\mathbf{h}_{\mathbf{E}}}\leq0.
\end{equation}

If $\nu=0$, all expressions to be estimated are zero. If
$\nu\neq0$, the projection property implies $q\in\partial K_x$, and
$\nu$ is an outward supporting normal at $q$. If $L$ is a
Lipschitz constant for $\Phi$ on the compact set traversed by $\mathbf{S}$ and its
projection, then~\eqref{eq:condicion-tangencia-convexa-parabolica} implies
\begin{align}
 \langle\nu,\Phi(x,t,\mathbf{S})\rangle_{\mathbf{h}_{\mathbf{E}}}
 &=\langle\nu,\Phi(x,t,q)\rangle_{\mathbf{h}_{\mathbf{E}}}
   +\langle\nu,\Phi(x,t,\mathbf{S})-\Phi(x,t,q)\rangle_{\mathbf{h}_{\mathbf{E}}}
 \nonumber\\
 &\leq L|\nu|_{\mathbf{h}_{\mathbf{E}}}^2=2L\rho(x,t).
 \label{eq:reaccion-apoyo-convexo}
\end{align}
Identities~\eqref{eq:derivada-temporal-distancia-convexo} and
\eqref{eq:sistema-principio-maximo-tensorial}, together with
\eqref{eq:parte-espacial-apoyo-convexo} and
\eqref{eq:reaccion-apoyo-convexo}, show that
\begin{equation}
\label{eq:derivada-maximo-distancia-convexo}
 \partial_t\rho(x,t)\leq2L\rho(x,t)
\end{equation}
at every spatial maximum of $\rho(\cdot,t)$.

It remains to pass from the spatial maximum to a function of one variable. Define
$\displaystyle Q(t)=\displaystyle\max_{x\in M}\rho(x,t)$. Its upper Dini derivative satisfies
\begin{equation}
\label{eq:dini-maximo-distancia-convexo}
 D^+Q(t)\leq
 \max_{\substack{x\in M\\\rho(x,t)=Q(t)}}\partial_t\rho(x,t)\leq2LQ(t).
\end{equation}
For the first inequality, choose $h_j\to 0^{+}$ and points $x_j$ where
$\rho(\cdot,t+h_j)$ attains its maximum. A subsequence converges to a maximum
of $\rho(\cdot,t)$, and
\[
 \frac{Q(t+h_j)-Q(t)}{h_j}
 \leq\frac{\rho(x_j,t+h_j)-\rho(x_j,t)}{h_j};
\]
uniform continuity of $\partial_t\rho$ near the positive time $t$
allows passage to the limit. The following Grönwall-type estimate for
upper Dini derivatives is proved by applying the first-maximum argument to
$e^{-2Lt}Q(t)$; thus,
\[
 Q(t)\leq e^{2L(t-\tau)}Q(\tau),
 \qquad 0<\tau<t\leq T.
\]
Since $Q$ is continuous and $Q(0)=0$, let $\tau\to 0^{+}$ to obtain
$Q(t)=0$. Therefore, $\mathbf{S}(x,t)\in K_x$ for every $(x,t)$.
\end{proof}

\begin{remark}
The proof remains valid if the connection depends smoothly on
$t$, provided that it is compatible with the same fixed bundle metric and
each connection preserves $\mathcal K$ under parallel transport. Indeed,
the spatial argument is performed at a fixed time, and the function
$\displaystyle\frac{1}{2}\operatorname{dist}(\mathbf{S},K_x)^2$ does not introduce a varying
metric upon time differentiation. If the bundle metric also varies,
one must first transport the problem to a bundle with a fixed metric or
include its derivative in the time identity.
\end{remark}

\begin{corollary}[Cones of nonnegative forms]
\label{cor:cono-formas-principio-maximo-tensorial}
Let $(M,\mathbf{g}_0)$ be a closed Riemannian manifold,
let $\mathbf{V}\longrightarrow M$ be a real vector bundle with a bundle
metric and compatible connection, and let
$\mathbf{S}(t)$ be a section of $S^2\mathbf{V}^*$ satisfying
\[
 \partial_t\mathbf{S}=\sum_{i,j=1}^m a^{ij}(\nabla^2\mathbf{S})_{ij}+\Phi(\mathbf{S}).
\]
Suppose that, for every form $Q\geq0$ and every $v\in \mathbf{V}_x$ with
$Q(v,v)=0$, we have
\[
 \Phi(Q)(v,v)\geq0.
\]
Also assume the regularity and uniform positivity of the
coefficients and the local Lipschitz condition on $\Phi$ from
Theorem~\ref{teo:principio-maximo-parabolico-tensorial}.
If $\mathbf{S}(0)\geq0$, then $\mathbf{S}(t)\geq0$ for as long as the solution exists.
\end{corollary}

\begin{proof}
The cone of positive semidefinite forms is closed, convex, and invariant
under parallel transport. We justify its supporting normals at a form \(Q\geq0\).
Identify symmetric forms with self-adjoint endomorphisms using the
metric on \(\mathbf V_x\). An outward normal functional is
represented by a self-adjoint matrix \(V\) satisfying, for every
\(Z\geq0\),
\[
 \operatorname{tr}(V(Z-Q))\leq0.
\]
Taking \(Z=0\) and \(Z=2Q\) gives
\(\operatorname{tr}(VQ)=0\). Taking
\(Z=t\,v\otimes v\), with arbitrary \(t\geq0\), yields
\(\langle Vv,v\rangle\leq0\); thus \(-V\geq0\).
By the finite-dimensional spectral theorem there is an orthonormal
basis \((e_i)\) in which
\[
 -V=\sum_i\mu_i e_i\otimes e_i,\qquad \mu_i\geq0,
 \qquad 0=-\operatorname{tr}(VQ)=\sum_i\mu_i Q(e_i,e_i).
\]
The summands are nonnegative. If \(\mu_i>0\), then
\(Q(e_i,e_i)=0\), and positive semidefiniteness implies
\(e_i\in\ker Q\). The normal functional is therefore
\[
 Z\longmapsto-\sum_{\mu_i>0}\mu_i Z(e_i,e_i).
\]
The stated condition gives for this functional
\(\operatorname{tr}(V\Phi(Q))
   =-\sum_{\mu_i>0}\mu_i\Phi(Q)(e_i,e_i)\leq0\).

The hypothesis is precisely
\eqref{eq:condicion-tangencia-convexa-parabolica}; apply
Theorem~\ref{teo:principio-maximo-parabolico-tensorial}.
\end{proof}

For classical treatments of parabolic theory with less regular
coefficients, see \cite[Chapter~7]{Evans} and
\cite[Chapter~3]{chow2004ricci}.

\chapter{Ricci flow}
\label{cap:flujo-ricci}
\index{Ricci flow}
\index{DeTurck trick}

A Riemannian metric determines how lengths, angles,
and curvature are measured. Ricci flow allows this rule of measurement to change
with time according to an equation built from its own Ricci
tensor. Rather than searching directly for a special metric, one starts from an
initial metric and studies which geometric features improve, which are
preserved, and how curvature may concentrate.

The analogy with heat is useful, but should not be taken literally. The
unknown is the metric itself and the equation is invariant under
diffeomorphisms. This symmetry makes its symbol degenerate in the
directions tangent to the orbits of the diffeomorphism action. The
DeTurck procedure fixes a gauge using a vector field and
turns the problem into a strictly parabolic system; a family of
diffeomorphisms then recovers the solution of the original flow
\cite{Hamilton1982,DeTurck1983}.

We begin with variation formulas for the metric, connection, and
curvature. We then identify the symbol degeneracy,
construct the Ricci--DeTurck system, and prove short-time existence and
uniqueness. Curvature evolution equations and
maximum principles lead to the first a priori estimates. Finally,
we study rescalings, normalized flow, Einstein metrics, and
solitons. These last solutions provide self-similar models
that will reappear when regions in which singularities form are magnified.

\begin{semblanzaHistorica}{Richard Hamilton and a new strategy for geometry}
Hamilton introduced Ricci flow in 1982 and proposed studying the topology
of a manifold by letting its geometry evolve. The strategy required
combining three levels: a parabolic equation to construct and regularize the
metric, Riemannian estimates to control curvature, and topological
arguments to interpret possible limits. The compactness and singularity
analysis of the next chapter continues this program.
\end{semblanzaHistorica}

Let $M$ be a closed smooth manifold of dimension $m$, with $m\geq2$. Whenever an assertion
requires connectedness, this will be stated explicitly. For
a tensor $\mathbf{T}$ write
$\Delta \mathbf{T}=\operatorname{tr}_{\mathbf{g}}(\nabla^2\mathbf{T})=-\nabla^*\nabla \mathbf{T}$. All time-dependent
norms, contractions, and covariant derivatives are
calculated with $\mathbf{g}(t)$ unless otherwise stated. Coordinate
indices of $M$ belong to $\{1,\ldots,m\}$. For higher-order derivatives,
retain the order in
Definition~\ref{derivada covariante de orden superior para haces vectoriales}:
\[
 (\nabla^2\mathbf{T})(X,Y,\ldots)
 =\bigl(\nabla_X(\nabla_Y\mathbf{T})
          -\nabla_{\nabla_XY}\mathbf{T}\bigr)(\ldots).
\]
Thus, in $(\nabla^2\mathbf{T})_{ab i_1\ldots i_q}$, the index $a$
corresponds to the outer derivative and $b$ to the inner one. Permutations of
$a$ and $b$ will be performed using the curvature identity of the cited
chapter. We also retain
$R_{ijk}{}^\ell=(R(\partial_i,\partial_j)\partial_k)^\ell$,
$\displaystyle \operatorname{Ric}_{ij}=\displaystyle\sum_{k=1}^mR_{kij}{}^k$, and
$\displaystyle R_{ijk\ell}=\displaystyle\sum_{q=1}^m g_{\ell q}R_{ijk}{}^q$. For a symmetric covariant
tensor $\mathbf{h}$ set
\[
(\operatorname{div}\mathbf{h})_i:=\sum_{j=1}^m\nabla^jh_{ij},
\qquad
\operatorname{div}^2\mathbf{h}:=\sum_{i=1}^m\nabla^i(\operatorname{div}\mathbf{h})_i.
\]
Denote by $\lambda_{\mathbf{g}}$ the Riemann--Lebesgue measure
determined by $\mathbf{g}$. In a chart, its weight relative to Lebesgue
measure is $\sqrt{\det(\mathbf{g})}$; thus no orientation of $M$ is required.
When $M$ is orientable and an orientation is fixed, integration with respect to
$\lambda_{\mathbf{g}}$ agrees with integration of the Riemannian volume
form. Identities containing
$\partial_t d\lambda_{\mathbf{g}(t)}$ are understood as identities of
measures depending smoothly on $t$, verified through these coordinate
weights.

\section{Variations of the metric and the Einstein--Hilbert functional}
\label{sec:variaciones-metricas-ricci}
\index{Einstein--Hilbert functional}

The following variation identities fix the signs that
will appear later. Let $\mathbf{g}(t)$ be a smooth family of metrics,
$\mathbf{g}(0)=\mathbf{g}$, and
$\mathbf{h}=\left.\partial_t\mathbf{g}(t)\right|_{t=0}$.

\begin{lemma}[Metric variation identities]
\label{lem:identidades-variacion-metrica}
Let $M$ be a smooth manifold of dimension $m$ with or without boundary and let
$\mathbf{g}(t)$ be a smooth family of metrics on $M$, with
$\mathbf{g}(0)=\mathbf{g}$ and
$\mathbf{h}=\left.\partial_t\mathbf{g}(t)\right|_{t=0}$.
Then
\begin{align}
 \left.\partial_tg^{ij}(t)\right|_{t=0}
 &=-h^{ij},
 \label{eq:variacion-inversa-metrica}\\
 \left.\partial_t d\lambda_{\mathbf{g}(t)}\right|_{t=0}
 &=\frac{1}{2}\operatorname{tr}_{\mathbf{g}}(\mathbf{h})\,d\lambda_{\mathbf{g}},
 \label{eq:variacion-volumen-metrica}\\
 \left.\partial_t\Gamma(\mathbf{g}(t))^k_{ij}\right|_{t=0}
 &=\frac{1}{2}\sum_{\ell=1}^m g^{k\ell}
 \bigl(\nabla_i h_{j\ell}+\nabla_j h_{i\ell}
       -\nabla_\ell h_{ij}\bigr),
 \label{eq:variacion-christoffel}\\
 \left.\partial_tR(\mathbf{g}(t))_{ijr}{}^s\right|_{t=0}
 &=\nabla_i\dot\Gamma^s_{jr}-\nabla_j\dot\Gamma^s_{ir},
 \label{eq:variacion-riemann-conexion}\\
 \left.\partial_t\operatorname{Ric}(\mathbf{g}(t))_{ij}\right|_{t=0}
 &=\sum_{k=1}^m\bigl(\nabla_k\dot\Gamma^k_{ij}-\nabla_i\dot\Gamma^k_{kj}\bigr),
 \label{eq:variacion-ricci-conexion}\\
 DR_{\mathbf{g}}[\mathbf{h}]
 &=-\langle\operatorname{Ric}_{\mathbf{g}},\mathbf{h}\rangle_{\mathbf{g}}
   +\operatorname{div}^2\mathbf{h}-\Delta_{\mathbf{g}}\operatorname{tr}_{\mathbf{g}} \mathbf{h}.
 \label{eq:variacion-curvatura-escalar}
\end{align}
Here $\dot\Gamma^k_{ij}$ denotes the expression in the third line.
\end{lemma}

\begin{proof}
Differentiating $\displaystyle \sum_{k=1}^m g^{ik}(t)g_{kj}(t)=\delta^i_j$ gives
\[
 \sum_{k=1}^m\left(
 \left.\partial_tg^{ik}(t)\right|_{t=0}g_{kj}+g^{ik}h_{kj}\right)=0,
\]
and multiplication by $g^{j\ell}$ proves
\eqref{eq:variacion-inversa-metrica}. In coordinates, Jacobi's identity
for the determinant gives
\[
 \left.\partial_t\sqrt{\det \mathbf{g}(t)}\right|_{t=0}
 =\frac{1}{2}\sqrt{\det \mathbf{g}}\,\sum_{i,j=1}^m g^{ij}h_{ij},
\]
whence \eqref{eq:variacion-volumen-metrica}.

If $\nabla^t$ is the Levi--Civita connection of $\mathbf{g}(t)$ and
\[
 \mathbf{A}(\mathbf{X},\mathbf{Y})=
 \left.\frac{d}{dt}\right|_{t=0}\nabla^t_{\mathbf{X}}\mathbf{Y},
\],
differentiating the Koszul formula and canceling the
terms containing $\nabla_XY$ gives
\[
 2\mathbf{g}(\mathbf{A}(\mathbf{X},\mathbf{Y}),\mathbf{Z})
 =(\nabla_{\mathbf{X}}\mathbf{h})(\mathbf{Y},\mathbf{Z})
 +(\nabla_{\mathbf{Y}}\mathbf{h})(\mathbf{X},\mathbf{Z})
 -(\nabla_{\mathbf{Z}}\mathbf{h})(\mathbf{X},\mathbf{Y}).
\]
In a coordinate frame, this is \eqref{eq:variacion-christoffel}. Differentiating
the definition of curvature and using that $\nabla$ is torsion-free yields
\[
 \left.\partial_t\mathbf{R}^t(\mathbf{X},\mathbf{Y})\mathbf{Z}\right|_{t=0}
 =(\nabla_{\mathbf{X}}\mathbf{A})(\mathbf{Y},\mathbf{Z})
 -(\nabla_{\mathbf{Y}}\mathbf{A})(\mathbf{X},\mathbf{Z}),
\]
which proves \eqref{eq:variacion-riemann-conexion}; contraction of the
first covariant index with the contravariant index gives
\eqref{eq:variacion-ricci-conexion}.

The contractions of \eqref{eq:variacion-christoffel} are
\[
 \sum_{k=1}^m\dot\Gamma^k_{kj}=\frac{1}{2}\nabla_j\operatorname{tr}_{\mathbf{g}} \mathbf{h},
 \qquad
 \sum_{i,j=1}^m g^{ij}\dot\Gamma^k_{ij}
 =\sum_{i=1}^m\nabla^ih_i{}^k-\frac{1}{2}\nabla^k\operatorname{tr}_{\mathbf{g}} \mathbf{h}.
\]
Therefore,
\[
 \sum_{i,j=1}^m g^{ij}\left.\partial_t\operatorname{Ric}(\mathbf{g}(t))_{ij}\right|_{t=0}
 =\operatorname{div}^2\mathbf{h}-\Delta_{\mathbf{g}}\operatorname{tr}_{\mathbf{g}} \mathbf{h}.
\]
Finally, differentiate $\displaystyle R(\mathbf{g}(t))=\displaystyle\sum_{i,j=1}^m g^{ij}(t)\operatorname{Ric}(\mathbf{g}(t))_{ij}$ and
use \eqref{eq:variacion-inversa-metrica}. This gives
\eqref{eq:variacion-curvatura-escalar}.
\end{proof}

For $s>\frac{m}{2}+2$, the set $\operatorname{Met}^s(M)$ of
$H^s(M,S^2T^*M)$ metrics is open in $H^s(M,S^2T^*M)$. Indeed, fixing a smooth metric
$\mathbf{g}_*$, the embedding $H^s(M,S^2T^*M)\hookrightarrow C^0(M,S^2T^*M)$ makes the function
\[
 \mathbf{g}\longmapsto
 \min_{x\in M}\min_{\substack{v\in T_xM\\|v|_{\mathbf{g}_*}=1}}
 \frac{\mathbf{g}_x(v,v)}{\mathbf{g}_{*,x}(v,v)}.
\]
continuous. Denote this minimum by $m(\mathbf g)$. Its existence is
justified in local frames of $\mathbf g_*$ over a finite cover
by compact sets: the unit sphere is identified with
a fixed Euclidean sphere and evaluation of $\mathbf g$ is
continuous. This is the fiberwise extrema argument of
Lemma~\ref{lema equivalencia metricas fibradas}.
Moreover, if $\mathbf g_1,\mathbf g_2$ are two symmetric fields
of class $H^s(M,S^2T^*M)$ and $|v|_{\mathbf g_*}=1$, Cauchy--Schwarz gives
\[
 |(\mathbf g_1-\mathbf g_2)_x(v,v)|
 \leq\|\mathbf g_1-\mathbf g_2\|_{\infty;\mathbf g_*}.
\]
Taking minima and then interchanging $\mathbf g_1$ and
$\mathbf g_2$ yields
\[
 |m(\mathbf g_1)-m(\mathbf g_2)|
 \leq\|\mathbf g_1-\mathbf g_2\|_{\infty;\mathbf g_*}
 \leq C\|\mathbf g_1-\mathbf g_2\|_{H^s(M,S^2T^*M)}.
\]
This proves continuity with respect to the Sobolev
norm. Since $M$ is compact, $m(\mathbf g)>0$ precisely
when $\mathbf g$ is positive definite in every fiber.
The condition $m(\mathbf g)>0$ therefore defines an open set. Inversion of positive matrices,
Sobolev products, and the local curvature formulas show that
\[
 \mathbf{g}\longmapsto \int_M R_{\mathbf{g}}\,d\lambda_{\mathbf{g}}
\]
is a smooth map from $\operatorname{Met}^s(M)$ to $\mathbb R$.

\begin{proposition}[First variation of Einstein--Hilbert]
\label{prop:primera-variacion-einstein-hilbert}
Let $M$ be a closed connected smooth manifold of dimension $m$ and let
$\mathbf{g}$ be a Riemannian metric on $M$.
The functional
\begin{equation}
 \mathcal S(\mathbf{g})=\int_M R_{\mathbf{g}}\,d\lambda_{\mathbf{g}}
 \label{eq:funcional-einstein-hilbert}
\end{equation}
satisfies
\begin{equation}
 D\mathcal S(\mathbf{g})[\mathbf{h}]
 =-\int_M\left\langle
 \operatorname{Ric}_{\mathbf{g}}-\frac{1}{2}R_{\mathbf{g}}\mathbf{g},\mathbf{h}
 \right\rangle_{\mathbf{g}}\,d\lambda_{\mathbf{g}}.
 \label{eq:primera-variacion-einstein-hilbert}
\end{equation}
If $m\geq3$, the critical points of $\mathcal S$ restricted to metrics
of fixed volume are precisely the Einstein metrics.
\end{proposition}

\begin{proof}
Equations \eqref{eq:variacion-volumen-metrica} and
\eqref{eq:variacion-curvatura-escalar} give
\[
 D\mathcal S(\mathbf{g})[\mathbf{h}]
 =\int_M\left(
 -\langle\operatorname{Ric},\mathbf{h}\rangle_{\mathbf{g}}
 +\frac{1}{2}R\operatorname{tr}_{\mathbf{g}} \mathbf{h}
 +\sum_{i=1}^m\nabla_i\left(\sum_{j=1}^m\nabla_jh^{ij}-\nabla^i\operatorname{tr}_{\mathbf{g}} \mathbf{h}\right)
 \right)d\lambda_{\mathbf{g}}.
\]
The integral of the divergence over the closed manifold $M$ is zero,
proving \eqref{eq:primera-variacion-einstein-hilbert}.

On the other hand,
\[
 D\operatorname{Vol}_{\mathbf{g}}[\mathbf{h}]
 =\frac{1}{2}\int_M\operatorname{tr}_{\mathbf{g}} \mathbf{h}\,d\lambda_{\mathbf{g}}
\]
is not the zero functional, since its value at $\mathbf{h}=\mathbf{g}$ is
$\frac{m}{2}\operatorname{Vol}(M,\mathbf{g})$. If $\mathbf{g}$ is critical on the
volume level set, choose $\mathbf{h}_0$ with $D\operatorname{Vol}_{\mathbf{g}}[\mathbf{h}_0]=1$.
Every tensor $\mathbf{h}$
decomposes as
\[
 \mathbf{h}=\bigl(\mathbf{h}-D\operatorname{Vol}_{\mathbf{g}}[\mathbf{h}]\mathbf{h}_0\bigr)
   +D\operatorname{Vol}_{\mathbf{g}}[\mathbf{h}]\mathbf{h}_0.
\]
The first summand belongs to the kernel of $D\operatorname{Vol}_{\mathbf{g}}$, where
$D\mathcal S_{\mathbf{g}}$ vanishes. Therefore,
$D\mathcal S_{\mathbf{g}}=cD\operatorname{Vol}_{\mathbf{g}}$, with
$c=D\mathcal S_{\mathbf{g}}[\mathbf{h}_0]$. The formulas for the two differentials imply
\[
 \operatorname{Ric}-\frac{1}{2}R\mathbf{g}=c\mathbf{g}
\]
after changing the value of the constant $c$ if necessary. Taking the trace
of this equality when $m\geq3$
shows that $R$ is constant; substituting it back gives
$\operatorname{Ric}=\lambda \mathbf{g}$ with constant $\lambda$. The reverse
implication follows directly from
\eqref{eq:primera-variacion-einstein-hilbert}, because vectors tangent
to the volume level set satisfy
$\displaystyle \int_M\operatorname{tr}_{\mathbf{g}} \mathbf{h}\,d\lambda_{\mathbf{g}}=0$.
\end{proof}

In dimension two,
Theorem~\ref{teo:gauss-bonnet-superficies-apendice} makes
$\mathcal S$ constant on each topological component of the space of metrics;
this is why the last
characterization requires $m\geq3$. The first variation formula,
however, remains valid in every dimension.

\section{Diffeomorphism invariance and symbol degeneracy}
\label{sec:invariancia-difeomorfismos-ricci}

\begin{definition}[Ricci flow]
\label{def:flujo-ricci}
Let $M$ be a smooth manifold with or without boundary. A \emph{Ricci flow} on
$M$ is a family of metrics $\mathbf{g}(t)$ satisfying
\begin{equation}
 \partial_t\mathbf{g}(t)=-2\operatorname{Ric}_{\mathbf{g}(t)}.
 \label{eq:flujo-ricci}
\end{equation}
\end{definition}

The Levi--Civita connection is natural: if $\phi$ is a diffeomorphism,
$\nabla^{\phi^*\mathbf{g}}=\phi^*\nabla^{\mathbf{g}}$. Taking curvature and contracting
gives
\begin{equation}
 \operatorname{Ric}_{\phi^*\mathbf{g}}=\phi^*\operatorname{Ric}_{\mathbf{g}}.
 \label{eq:naturalidad-ricci}
\end{equation}
Thus pullback by a time-independent diffeomorphism takes
solutions to solutions. For a family $\phi_t$, if
\[
 \mathbf{X}_t=(\partial_t\phi_t)\circ\phi_t^{-1},
\],
the pullback differentiation rule is
\begin{equation}
 \frac{d}{dt}(\phi_t^*\mathbf{g}(t))
 =\phi_t^*(\partial_t\mathbf{g}(t)+\mathcal L_{\mathbf{X}_t}\mathbf{g}(t)).
 \label{eq:derivada-pullback-tiempo}
\end{equation}
It is proved first for functions and one-forms, then using the
Leibniz rule for tensor products.

Naturality explains why the right-hand side of
\eqref{eq:flujo-ricci} is not strictly parabolic. With
$\displaystyle (\delta_{\mathbf{g}}\mathbf{h})_j=-\displaystyle\sum_{i=1}^m\nabla^ih_{ij}$, define the formal adjoint of
$\delta_{\mathbf{g}}$
by
\begin{equation}
 (\delta_{\mathbf{g}}^*\boldsymbol{\omega})_{ij}
 =\frac{1}{2}\bigl(\nabla_i\omega_j+\nabla_j\omega_i\bigr).
 \label{eq:definicion-delta-estrella-ricci}
\end{equation}
If $\nabla^*\nabla=-\operatorname{tr}_{\mathbf{g}}(\nabla^2)$,
Lemma~\ref{lem:identidades-variacion-metrica}, followed by commutation of
covariant derivatives, gives
\begin{equation}
 2D\operatorname{Ric}_{\mathbf{g}}[\mathbf{h}]
 =\nabla^*\nabla \mathbf{h}-\nabla^2\operatorname{tr}_{\mathbf{g}} \mathbf{h}
  -2\delta_{\mathbf{g}}^*\delta_{\mathbf{g}}\mathbf{h}
  +\mathcal Q_{\mathbf{g}}(\mathbf{h}),
 \label{eq:linealizacion-ricci}
\end{equation}
where the zeroth-order term is
\begin{equation}
 \mathcal Q_{\mathbf{g}}(\mathbf{h})_{ij}
 =\sum_{k=1}^m\operatorname{Ric}_i{}^k h_{kj}
  +\sum_{k=1}^m\operatorname{Ric}_j{}^k h_{ki}
  -2\sum_{k,\ell=1}^mR_{kij\ell}h^{k\ell}.
 \label{eq:termino-cero-linealizacion-ricci}
\end{equation}
Before commuting, substitution of
\eqref{eq:variacion-christoffel} into
\eqref{eq:variacion-ricci-conexion} gives the tensor identity
\begin{equation}
 2D\operatorname{Ric}_{\mathbf{g}}[\mathbf{h}]_{ij}
 =\sum_{k,a=1}^m g^{ka}\bigl((\nabla^2\mathbf{h})_{a i j k}+(\nabla^2\mathbf{h})_{a j i k}\bigr)
  +(\nabla^*\nabla \mathbf{h})_{ij}
  -(\nabla^2\operatorname{tr}_{\mathbf{g}} \mathbf{h})_{ij}.
 \label{eq:linealizacion-ricci-indices}
\end{equation}
To see where $\mathcal Q_{\mathbf{g}}$ comes from, apply the
commutation identity to the two slots of $\mathbf{h}$:
\begin{align*}
 &\sum_{a,k=1}^m g^{ak}
 \bigl((\nabla^2\mathbf{h})_{a i j k}
              -(\nabla^2\mathbf{h})_{i a j k}\bigr)\\
 &\quad=-\sum_{a,k,p=1}^m g^{ak}
       \bigl(R_{aij}{}^p h_{pk}+R_{aik}{}^p h_{jp}\bigr)\\
 &\quad=\sum_{p=1}^m\operatorname{Ric}_i{}^p h_{pj}
              -\sum_{k,\ell=1}^mR_{kij\ell}h^{k\ell}.
\end{align*}
The same equality with $i$ and $j$ interchanged produces the other Ricci
term. The Riemann symmetries and symmetry of $\mathbf{h}$ make
the two remaining contractions equal. The derivatives already in the order
$i,a$ and $j,a$ form $-2\delta_{\mathbf{g}}^*\delta_{\mathbf{g}}\mathbf{h}$.
This proves \eqref{eq:linealizacion-ricci} and
\eqref{eq:termino-cero-linealizacion-ricci} with the chosen differentiation order.

If $\mathbf{h}=\mathcal L_{\mathbf{X}}\mathbf{g}$, its symbol part is
$h_{ij}=\xi_iX_j+\xi_jX_i$, and substitution into
\eqref{eq:linealizacion-ricci-indices} annihilates the symbol of
$D\operatorname{Ric}$. This is precisely the infinitesimal direction of the
diffeomorphism orbit. Thus strongly parabolic theory cannot be applied
directly to \eqref{eq:flujo-ricci}.

\section{The Ricci--DeTurck system}
\label{sec:sistema-ricci-deturck}

The preceding degeneracy reflects the freedom to change coordinates during
evolution rather than a loss of diffusion. The DeTurck procedure
removes this freedom by comparing the connection of the unknown metric with
a fixed connection. Fix a smooth metric $\bar{\mathbf{g}}$, denote its Levi--Civita connection by $\bar\nabla$,
and write
\[
 \mathbf{A}(\mathbf{g})^k_{ij}=\Gamma(\mathbf{g})^k_{ij}-\Gamma(\bar{\mathbf{g}})^k_{ij},
 \qquad
 \mathbf{W}_{\bar{\mathbf{g}}}(\mathbf{g})^k=\sum_{i,j=1}^m g^{ij}\mathbf{A}(\mathbf{g})^k_{ij}.
\]
The difference of connections is tensorial, so
$\mathbf{W}_{\bar{\mathbf{g}}}(\mathbf{g})$ is
a globally defined vector field. The modified system is
\begin{equation}
 \partial_t\mathbf{g}=-2\operatorname{Ric}_{\mathbf{g}}
 +\mathcal L_{\mathbf{W}_{\bar{\mathbf{g}}}(\mathbf{g})}\mathbf{g}.
 \label{eq:ricci-deturck}
\end{equation}

\begin{proposition}[Parabolicity of the DeTurck system]
\label{prop:parabolicidad-ricci-deturck}
Let $M$ be a smooth manifold of dimension $m$ with or without boundary, let
$\bar{\mathbf{g}}$ be a fixed Riemannian metric, and consider the Ricci--DeTurck
system \eqref{eq:ricci-deturck}.
The right-hand side of \eqref{eq:ricci-deturck} is a quasilinear operator
whose principal part is locally strongly parabolic. More precisely,
relative to $\bar\nabla$ it can be written as
\begin{equation}
 -2\operatorname{Ric}(\mathbf{g})_{ij}
 +(\mathcal L_{\mathbf{W}_{\bar{\mathbf{g}}}(\mathbf{g})}\mathbf{g})_{ij}
 =\sum_{a,b=1}^m g^{ab}(\bar\nabla^2\mathbf{g})_{abij}
 +P_{ij}(\mathbf{g},\mathbf{g}^{-1},\bar\nabla \mathbf{g},\overline{\operatorname{Rm}}),
 \label{eq:forma-local-ricci-deturck}
\end{equation}
where $P$ is smooth, contains no derivatives of $\mathbf{g}$ of order greater than one, and
is polynomial in $\bar\nabla \mathbf{g}$. The Fourier symbol of the linearization
at a covariable $\xi\neq0$ is
\begin{equation}
 -\sum_{a,b=1}^m g^{ab}\xi_a\xi_b\operatorname{Id}_{S^2T^*M}.
 \label{eq:simbolo-ricci-deturck}
\end{equation}
On each compact subset of $M$, the parabolicity bound can be chosen uniformly
for metrics sufficiently close in $C^0$ to a fixed metric.
In particular, if $M$ is compact, this bound is global.
\end{proposition}

\begin{proof}
The difference-of-connections formula is
\[
 \mathbf{A}(\mathbf{g})^k_{ij}
 =\frac{1}{2}\sum_{\ell=1}^m g^{k\ell}
 \bigl(\bar\nabla_i g_{j\ell}+\bar\nabla_jg_{i\ell}
       -\bar\nabla_\ell g_{ij}\bigr).
\]
Expressing the curvature of $\mathbf{g}$ using $\bar\nabla$ gives
\[
 \mathbf{R}(\mathbf{g})_{ijk}{}^\ell
 =\overline R_{ijk}{}^\ell
 +\bar\nabla_i\mathbf{A}(\mathbf{g})^\ell_{jk}
 -\bar\nabla_j\mathbf{A}(\mathbf{g})^\ell_{ik}
 +\sum_{p=1}^m\bigl(
 \mathbf{A}(\mathbf{g})^p_{jk}\mathbf{A}(\mathbf{g})^\ell_{ip}
 -\mathbf{A}(\mathbf{g})^p_{ik}\mathbf{A}(\mathbf{g})^\ell_{jp}\bigr).
\]
The contraction defining Ricci contains second derivatives of $\mathbf{g}$. On the other
hand,
\[
 (\mathcal L_{\mathbf{W}}\mathbf{g})_{ij}=\nabla_iW_j+\nabla_jW_i.
\]
At a point where normal coordinates for $\bar{\mathbf{g}}$ are chosen, the
parts containing two derivatives of $\mathbf{g}$ are
\begin{align*}
 (-2\operatorname{Ric}(\mathbf{g}))_{ij}^{(2)}
 =\sum_{a,b=1}^m g^{ab}\bigl(&
 (\bar\nabla^2\mathbf{g})_{abij}
 +(\bar\nabla^2\mathbf{g})_{ijab}
 -(\bar\nabla^2\mathbf{g})_{aibj}
 -(\bar\nabla^2\mathbf{g})_{ajbi}\bigr),\\
 (\mathcal L_{\mathbf{W}}\mathbf{g})_{ij}^{(2)}
 =\sum_{a,b=1}^m g^{ab}\bigl(&
 (\bar\nabla^2\mathbf{g})_{iabj}
 +(\bar\nabla^2\mathbf{g})_{jabi}
 -(\bar\nabla^2\mathbf{g})_{ijab}\bigr).
\end{align*}
The three gauge terms cancel in pairs, leaving as the only second-order
term $\displaystyle \sum_{a,b=1}^m g^{ab}(\bar\nabla^2\mathbf{g})_{abij}$. Commutations of
$\bar\nabla$ produce $\overline{\operatorname{Rm}}*\mathbf{g}$, and all other
terms are contractions of $\mathbf{g}$, $\mathbf{g}^{-1}$, $\bar\nabla \mathbf{g}$, and
$\overline{\operatorname{Rm}}$. This proves
\eqref{eq:forma-local-ricci-deturck}.

To apply the divergence-form theory, observe the exact identity
\begin{equation}
 \sum_{a,b=1}^m g^{ab}(\bar\nabla^2\mathbf{g})_{ab}
 =\sum_{a,b=1}^m\bar\nabla_a\bigl(g^{ab}\bar\nabla_b \mathbf{g}\bigr)
  -\sum_{a,b=1}^m(\bar\nabla_a g^{ab})\bar\nabla_b \mathbf{g}.
 \label{eq:deturck-forma-divergencia}
\end{equation}
The last summand contains only first derivatives of $\mathbf{g}$. The principal
coefficient is therefore
\[
 A^{ab}(\mathbf{g})=g^{ab}\operatorname{Id}_{S^2T^*M}.
\]
As an endomorphism of $T^*M\otimes S^2T^*M$ measured with
$\bar{\mathbf{g}}$, it reads
$\displaystyle A_a{}^b=\displaystyle\sum_{c=1}^m\bar g_{ac}g^{cb}\operatorname{Id}_{S^2T^*M}$;
raising the first index with $\bar{\mathbf{g}}$ recovers $A^{ab}$.
This coefficient satisfies $A^{ab*}=A^{ba}$. To specify uniformity,
fix a compact set $K\subseteq M$ and a continuous metric $\mathbf{g}_*$.
The comparison on compact sets proved after
Lemma~\ref{lema cg}, whose argument also applies to continuous
metrics, gives $0<c_*\leq C_*<\infty$ such that
\[
 c_*\bar{\mathbf g}_x(v,v)
 \leq(\mathbf g_*)_x(v,v)
 \leq C_*\bar{\mathbf g}_x(v,v)
 \qquad(x\in K,\ v\in T_xM).
\]
If $\displaystyle\|\mathbf g-\mathbf g_*\|_{C^0(K,S^2T^*M;\bar{\mathbf g})}<c_*/2$,
Cauchy--Schwarz for the tensor inner
product gives, for $x\in K$ and $v\in T_xM$,
\[
 |(\mathbf g-\mathbf g_*)_x(v,v)|
 \leq|\mathbf g-\mathbf g_*|_{\bar{\mathbf g}}(x)
      |v|_{\bar{\mathbf g}}^2
 <\frac{c_*}{2}\bar{\mathbf g}_x(v,v)
 \qquad(v\neq0).
\]
Indeed, in a $\bar{\mathbf g}$--orthonormal basis, apply
Cauchy--Schwarz to
$\displaystyle\sum_{a,b=1}^{m}(g-g_*)_{ab}v^av^b$ and use
$\displaystyle\sum_{a,b=1}^{m}(v^av^b)^2=(\displaystyle\sum_{a=1}^{m}(v^a)^2)^2$.
Adding this bound to the comparison for $\mathbf g_*$ gives
\[
 \frac{c_*}{2}\bar{\mathbf{g}}\leq\mathbf{g}
 \leq\left(C_*+\frac{c_*}{2}\right)\bar{\mathbf{g}}
 \qquad\text{on }K.
\]
The dual comparison in
Proposition~\ref{prop:comparacion-metricas-duales-tensores-volumen}
gives, with $\lambda=(C_*+c_*/2)^{-1}>0$,
\[
 \mathbf g_x^{-1}(\xi,\xi)
 \geq\lambda\bar{\mathbf g}_x^{-1}(\xi,\xi)
 \qquad(x\in K,\ \xi\in T_x^*M).
\]
To apply this inequality to a covector with values in
$S^2T_x^*M$, choose a basis $F_1,\ldots,F_d$ of
that space orthonormal for the metric induced by $\bar{\mathbf g}$,
where $d=m(m+1)/2$.
Write $\zeta_a=\displaystyle\sum_{\alpha=1}^d z_a^\alpha F_\alpha$
in a $\bar{\mathbf g}$--orthonormal tangent frame. For each
$\alpha$, the numbers $z_a^\alpha$ are the components of a
covector to which the preceding bound applies. Thus
\[
 \begin{aligned}
 \sum_{a,b=1}^m
 \left\langle A^{ab}(\mathbf g)\zeta_b,\zeta_a\right\rangle_{\bar{\mathbf g}}
 &=\sum_{\alpha=1}^d\sum_{a,b=1}^m
   g^{ab}z_b^\alpha z_a^\alpha\\
 &\geq\lambda\sum_{\alpha=1}^d\sum_{a=1}^m|z_a^\alpha|^2
 =\lambda|\zeta|_{\bar{\mathbf g}}^2.
 \end{aligned}
\]
If $M$ is compact, take $K=M$. Sobolev embedding into $C^0$ allows
the neighborhood in $\operatorname{Met}^s(M)$ to be chosen so that the
preceding bound holds; in particular, this applies for $s>m/2+2$, as in this
chapter. For closed manifolds, this gives the uniform Legendre
condition of Theorem~\ref{teo:existencia-parabolica-cuaslineal};
the permitted open subset of bundle values is the cone of positive definite
bilinear forms. If $M$ is noncompact, a comparison
$\mathbf{g}\leq C\bar{\mathbf{g}}$ on all of $M$, with $C$ independent of the
point, likewise gives a lower bound for the principal part
$\lambda=C^{-1}$. Without uniform comparison, pointwise positivity alone does not imply a
global constant.

Linearization of the principal term in the direction $\mathbf{h}$ is
$\displaystyle\sum_{a,b=1}^m g^{ab}(\bar\nabla^2\mathbf{h})_{abij}$ plus terms of order at most one
in $\mathbf{h}$. The substitution $\bar\nabla_a\mapsto i\xi_a$ gives
\eqref{eq:simbolo-ricci-deturck}. Since $\displaystyle \sum_{a,b=1}^m g^{ab}\xi_a\xi_b>0$, the system is
locally strongly parabolic, and the preceding bounds provide
uniformity in the compact case. The presence of a boundary does not alter
this symbol calculation; the theory of
Chapter~\ref{cap:teoria-parabolica-cerradas} is applied here with $M$ closed.
\end{proof}

The uniqueness argument will use harmonic map heat flow. To
construct it through an equation on a trivial bundle, we must check
that the ambient solution remains in the embedded submanifold. The following
estimate makes this invariance explicit.

\begin{lemma}[Remaining in a submanifold through tubular distance]
\label{lem:permanencia-subvariedad-flujo-armonico}
Let $N$ be a compact embedded submanifold without boundary of a Riemannian manifold
$(\mathcal N,\mathbf{h})$, and let $\mathcal T$ be a tubular neighborhood with
retraction $\pi\colon\mathcal T\to N$. Set
$\rho(y)=\frac12\operatorname{dist}_{\mathbf{h}}(y,N)^2$. Let
$\widetilde{\boldsymbol{\mathcal A}}$ be a smooth extension to $\mathcal T$ of
the second fundamental form $\boldsymbol{\mathcal A}$ of $N$, defined
as follows. For $y=\exp_z^\perp\nu$, with $z=\pi(y)$, let
$\mathsf P_{z,y}\colon T_z\mathcal N\to T_y\mathcal N$ be parallel
transport along the normal geodesic. Set
\[
 \widetilde{\boldsymbol{\mathcal A}}_y(V,W)
 :=\mathsf P_{z,y}\bigl(
       \boldsymbol{\mathcal A}_z(d\pi_yV,d\pi_yW)\bigr),
 \qquad V,W\in T_y\mathcal N.
\]
If a smooth map
$U\colon M\times[0,T]\to\mathcal T$ satisfies
\begin{equation}
 \partial_tU=\tau_{\mathbf{g}(t),\mathbf{h}}(U)
 -\sum_{i,j=1}^m g^{ij}(t)\widetilde{\boldsymbol{\mathcal A}}_U
       (\partial_iU,\partial_jU),
 \label{eq:sistema-ambiente-lema-permanencia}
\end{equation},
there exists a constant $C$, uniform while $U$ ranges over a compact subset of
$\mathcal T$ and $\mathbf{g}(t)$ stays in a uniformly equivalent
family, such that
\begin{equation}
 (\partial_t-\Delta_{\mathbf{g}(t)})(\rho\circ U)
 \leq C\bigl(1+|dU|_{\mathbf{g}(t),\mathbf{h}}^2\bigr)(\rho\circ U).
 \label{eq:distancia-flujo-armonico}
\end{equation}
In particular, if $M$ is closed and $U(\cdot,0)$ takes values in $N$, then
$U(M\times[0,T])\subseteq N$.
\end{lemma}

\begin{proof}
The chain rule for the tension field of a map gives
\begin{align*}
 (\partial_t-\Delta_{\mathbf{g}(t)})(\rho\circ U)
 ={}&-d\rho_U\!\left(\sum_{i,j=1}^m g^{ij}
 \widetilde{\boldsymbol{\mathcal A}}_U
 (\partial_iU,\partial_jU)\right)\\
 &-\sum_{i,j=1}^m g^{ij}\operatorname{Hess}_{\mathbf{h}}\rho_U
 (\partial_iU,\partial_jU).
\end{align*}
To estimate the right-hand side, use Fermi coordinates
$y=(z^a,\nu^\alpha)$, where $z=\pi(y)$ and $\nu$ is normal to $N$ at $z$.
Write $d=\dim N$ and $c=\dim\mathcal N-d$.
Choose a local orthonormal normal frame
$(\mathbf{n}_\alpha)_{\alpha=1}^c$ and set
$\displaystyle \nu=\displaystyle\sum_{\alpha=1}^c\nu^\alpha\mathbf{n}_\alpha(z)$.
In these coordinates, $a,b\in\{1,\ldots,d\}$ and
$\alpha,\beta\in\{1,\ldots,c\}$. On a compact subset
of the tubular neighborhood, the Taylor expansions are uniform and give
\begin{align}
 \rho(z,\nu)&=\frac12\sum_{\alpha=1}^c(\nu^\alpha)^2,
 &\partial_a\rho&=0,
 &\partial_\alpha\rho&=\nu^\alpha,
 \label{eq:expansion-rho-tubular}\\
 (\operatorname{Hess}\rho)_{ab}
 &=-\langle\nu,\boldsymbol{\mathcal A}_{ab}\rangle+O(|\nu|^2),
 &(\operatorname{Hess}\rho)_{a\beta}&=O(|\nu|),
 &(\operatorname{Hess}\rho)_{\alpha\beta}
 &=\delta_{\alpha\beta}+O(|\nu|).
 \label{eq:expansion-hessiana-rho-tubular}
\end{align}
The first line is exact: normal geodesics minimize distance
to $N$ in the tubular neighborhood. For the second line, use
$\displaystyle (\operatorname{Hess}\rho)_{AB}
=\partial_A\partial_B\rho-\displaystyle\sum_{C=1}^{d+c}\Gamma^C_{AB}\partial_C\rho$.
The indices $A,B$ range over $\{1,\ldots,d+c\}$. On the zero section, $\Gamma^\alpha_{ab}$ are the normal components of
$\boldsymbol{\mathcal A}_{ab}$; multiplying them by $\nu^\alpha$ gives the
linear term of the tangential block. The normal block starts with the
identity and the mixed block vanishes on $N$, justifying the remainder orders
by Taylor's formula.
Here $\boldsymbol{\mathcal A}_{ab}$ is the second fundamental form of
$N$ on the zero section. By the way this form was extended,
transport $V\in T_y\mathcal N$ back to $T_z\mathcal N$ and
decompose it as $V^\top+V^\perp$ according to $T_zN\oplus(T_zN)^\perp$.
Smoothness of $d\pi$ and its value on the zero section give
$d\pi_yV=V^\top+O(|\nu|)|V|$. Moreover,
$d\rho_y(\mathsf P_{z,y}Z)=\langle\nu,Z\rangle$ for
$Z\in T_z\mathcal N$ by the first variation of distance. Thus
\begin{equation}
 d\rho_y\!\left(
 \widetilde{\boldsymbol{\mathcal A}}_y(V,V)\right)
 =\langle\nu,\boldsymbol{\mathcal A}_z(V^\top,V^\top)\rangle
 +O(|\nu|^2)|V|_{\mathbf{h}}^2.
 \label{eq:expansion-segunda-forma-tubular}
\end{equation}

Fix a point of the domain and choose there an orthonormal frame
$(\mathbf{e}_i)_{i=1}^m$ for $\mathbf{g}(t)$. If
$V_i=dU(\mathbf{e}_i)=V_i^\top+V_i^\perp$, the two preceding formulas
show, for each $i$, that
\begin{align*}
 &-d\rho_U\!\left(
 \widetilde{\boldsymbol{\mathcal A}}_U(V_i,V_i)\right)
 -\operatorname{Hess}_{\mathbf{h}}\rho_U(V_i,V_i)\\
 &\qquad=-|V_i^\perp|_{\mathbf{h}}^2+\mathcal E_i,
\end{align*}
with
\[
 |\mathcal E_i|
 \leq C|\nu|\,|V_i^\perp|_{\mathbf{h}}|V_i|_{\mathbf{h}}
      +C|\nu|^2\bigl(1+|V_i|_{\mathbf{h}}^2\bigr).
\]
Indeed, the term
$\langle\nu,\boldsymbol{\mathcal A}(V_i^\top,V_i^\top)\rangle$ coming
from $d\rho(\widetilde{\boldsymbol{\mathcal A}})$ cancels with the
linear term of opposite sign in
$\operatorname{Hess}\rho(V_i^\top,V_i^\top)$. The normal block of the
Hessian leaves $-|V_i^\perp|^2$, while the mixed block and Taylor
remainders produce the stated errors. Summing over $i$ gives
\begin{align}
 &-d\rho_U\!\left(\sum_{i,j=1}^m g^{ij}
 \widetilde{\boldsymbol{\mathcal A}}_U
 (\partial_iU,\partial_jU)\right)
 -\sum_{i,j=1}^m g^{ij}\operatorname{Hess}_{\mathbf{h}}\rho_U
 (\partial_iU,\partial_jU)\nonumber\\
 &\qquad=-|(dU)^\perp|_{\mathbf{g}(t),\mathbf{h}}^2+\mathcal E,
 \label{eq:cancelacion-normal-flujo-armonico}
\end{align}
where
\begin{equation}
 |\mathcal E|
 \leq C|\nu|\,|(dU)^\perp|\,|dU|
       +C|\nu|^2\bigl(1+|dU|^2\bigr).
 \label{eq:error-tubular-flujo-armonico}
\end{equation}

Young's inequality absorbs the cross term:
\[
 C|\nu|\,|(dU)^\perp|\,|dU|
 \leq\frac12|(dU)^\perp|^2+C|\nu|^2|dU|^2.
\]
In these coordinates, $|\nu|^2=2\rho(z,\nu)$. The three preceding estimates
imply \eqref{eq:distancia-flujo-armonico}. On an interval
where $U$ is smooth, $|dU|$ is bounded because $M$ is compact. The
maximum principle applied to
$e^{-Ct}(\rho\circ U)$ with zero initial data gives
$\rho\circ U=0$.
\end{proof}

\begin{corollary}[Harmonic map heat flow with varying domain metric]
\label{cor:flujo-armonico-aplicaciones-global}
Let $M$ be a closed smooth manifold. Let $\mathbf{g}(t)$ be a smooth family of metrics on $M$, let $\bar{\mathbf{g}}$ be a fixed smooth
metric, and let $F_0\colon M\to M$ be a smooth map. There exists $T_0>0$ for
which the problem
\begin{equation}
 \partial_tF=\tau_{\mathbf{g}(t),\bar{\mathbf{g}}}(F),
 \qquad F(0)=F_0,
 \label{eq:flujo-armonico-global}
\end{equation}
has a unique smooth solution. The time can be chosen uniformly
when the domain metrics and their inverses satisfy uniform bounds
for ellipticity and regularity, and the data $\iota\circ F_0$ are bounded
in an $H^s(M,\mathbb R^N)$ norm, with $s>m/2+2$, for a fixed embedding $\iota$ of the target.
The smooth geometry of the target remains fixed.
\end{corollary}

\begin{proof}
Theorem~\ref{teo:whitney-encaje-apendice} gives an embedding
$\iota\colon M\hookrightarrow\mathbb R^N$. On a neighborhood of $\iota(M)$,
choose an ambient metric $\mathbf{h}$ making the inclusion
$\iota\colon(M,\bar{\mathbf{g}})\to(\mathbb R^N,\mathbf{h})$ isometric. To
construct it, on the normal bundle of a Euclidean embedding take the orthogonal
sum of $\bar{\mathbf{g}}$ in the horizontal directions and a
smooth metric in the normal directions. A fixed connection defines this splitting;
the tubular diffeomorphism transports the metric to a neighborhood of
$\iota(M)$. Its restriction to the zero section is the prescribed one. Let
$\mathcal T$ be a tubular neighborhood, let
$\pi\colon\mathcal T\to\iota(M)$ be its retraction, and let
$\boldsymbol{\mathcal A}$ be the second fundamental form of $\iota(M)$.
Shrinking $\mathcal T$, smoothly extend the map
\[
 (y,V,W)\longmapsto
 \mathsf P_{\pi(y),y}
 \bigl(\boldsymbol{\mathcal A}_{\pi(y)}(d\pi_yV,d\pi_yW)\bigr)
\]
to an open subset of $\mathbb R^N$ and denote the extension by
$\widetilde{\boldsymbol{\mathcal A}}$. Existence of the tubular neighborhood
and retraction was proved in
Theorem~\ref{teo:variedades-riemannianas-existencia-de-vecindades-tubulares}.

For $U_0=\iota\circ F_0$, solve on the trivial bundle
$M\times\mathbb R^N$ the system
\begin{equation}
 \partial_tU=\tau_{\mathbf{g}(t),\mathbf{h}}(U)
 -\sum_{i,j=1}^m g^{ij}(t)\widetilde{\boldsymbol{\mathcal A}}_U
       (\partial_iU,\partial_jU),
 \qquad U(0)=U_0.
 \label{eq:extension-ambiente-flujo-armonico}
\end{equation}
In ambient-space coordinates, the principal part is
$\displaystyle \sum_{i,j=1}^m g^{ij}(t)\partial_i\partial_jU^\alpha$ and the remaining terms depend
smoothly on $(x,t,U,dU)$.
Corollary~\ref{cor:existencia-parabolica-cuaslineal-tiempo} gives a unique
smooth solution on a short interval. Shrinking this interval, its image
stays in $\mathcal T$.

Apply
Lemma~\ref{lem:permanencia-subvariedad-flujo-armonico} with
$N=\iota(M)$. Since $U_0$ takes values in $N$, we have
$U(M\times[0,T_0])\subseteq\iota(M)$. Thus there exists a smooth
map $F$ such that $U=\iota\circ F$. The Gauss formula states that the
normal component of the tension field of $\iota\circ F$ is
$\displaystyle \sum_{i,j=1}^m g^{ij}\boldsymbol{\mathcal A}(dF(\partial_i),dF(\partial_j))$; restricting
\eqref{eq:extension-ambiente-flujo-armonico} to the submanifold
cancels this term, leaving exactly
\eqref{eq:flujo-armonico-global}.

Composing two intrinsic solutions with $\iota$ gives two solutions
of the same ambient system, so they agree. Finally, the existence constants
depend only on the uniform bounds for $\mathbf{g}$ and
$\mathbf{g}^{-1}$ and finitely many of their derivatives, the fixed geometry
of the tubular neighborhood, and an $H^s(M,\mathbb R^N)$ bound for $U_0=\iota\circ F_0$.
The embedding $H^s(M,\mathbb R^N)\hookrightarrow C^1(M,\mathbb R^N)$ also controls the values and first
derivatives occurring in the nonlinearity. This proves the uniform assertion.
\end{proof}

\begin{theorem}[Short-time existence and uniqueness]
\label{teo:existencia-unicidad-flujo-ricci}
For every smooth metric $\mathbf{g}_0$ on a closed smooth manifold $M$, there exist
$T>0$ and a unique smooth solution
\[
\mathbf{g}\colon M\times[0,T]\longrightarrow S^2T^*M
\]
of Ricci flow with $\mathbf{g}(0)=\mathbf{g}_0$. A smooth family of smooth initial
metrics gives, locally in the parameter, a smooth family
of solutions on a common interval. The interval can be chosen
uniformly over compact subsets of the parameter space. In particular, for $\tau>0$,
dependence is smooth in every $C^k$ norm on $M\times[\tau,T]$, as long as
the existence interval is common.
\end{theorem}

\begin{proof}
Take $\bar{\mathbf{g}}=\mathbf{g}_0$.
Proposition~\ref{prop:parabolicidad-ricci-deturck} and
Theorem~\ref{teo:existencia-parabolica-cuaslineal} give a unique
solution $\widetilde{\mathbf{g}}(t)$ of
\begin{equation}
 \partial_t\widetilde{\mathbf{g}}
 =-2\operatorname{Ric}_{\widetilde{\mathbf{g}}}
  +\mathcal L_{\mathbf{W}_{\bar{\mathbf{g}}}(\widetilde{\mathbf{g}})}
   \widetilde{\mathbf{g}},
 \qquad \widetilde{\mathbf{g}}(0)=\mathbf{g}_0.
 \label{eq:problema-ricci-deturck-inicial}
\end{equation}
Positive definiteness is an open condition; continuity in time
allows $T$ to be reduced so that $\widetilde{\mathbf{g}}(t)$ remains a metric.
Let $\phi_t$ be the nonautonomous flow determined by
\begin{equation}
 \partial_t\phi_t
 =-\mathbf{W}_{\bar{\mathbf{g}}}(\widetilde{\mathbf{g}}(t))\circ\phi_t,
 \qquad \phi_0=\operatorname{id}_M.
 \label{eq:flujo-difeomorfismos-deturck}
\end{equation}
The field is smooth on $M\times[0,T]$. Compactness of $M$ ensures that the
flow exists on the whole interval, and each $\phi_t$ is a diffeomorphism. With
$\mathbf{g}(t)=\phi_t^*\widetilde{\mathbf{g}}(t)$, equations
\eqref{eq:derivada-pullback-tiempo},
\eqref{eq:problema-ricci-deturck-inicial}, and
\eqref{eq:naturalidad-ricci} give
\begin{align*}
 \partial_t\mathbf{g}
 &=\phi_t^*\left(
 -2\operatorname{Ric}_{\widetilde{\mathbf{g}}}
 +\mathcal L_{\mathbf{W}_{\bar{\mathbf{g}}}(\widetilde{\mathbf{g}})}\widetilde{\mathbf{g}}
 -\mathcal L_{\mathbf{W}_{\bar{\mathbf{g}}}(\widetilde{\mathbf{g}})}\widetilde{\mathbf{g}}
 \right)\\
 &=-2\operatorname{Ric}_{\phi_t^*\widetilde{\mathbf{g}}}
 =-2\operatorname{Ric}_{\mathbf{g}}.
\end{align*}
This proves existence.

Uniqueness does not follow by applying uniqueness of the DeTurck system to two
Ricci flows: the gauge-fixing diffeomorphisms depend on the
solution. To address this difficulty, consider two solutions
$\mathbf{g}_1(t)$ and $\mathbf{g}_2(t)$ with the same initial data and solve, for $a=1,2$,
the harmonic map heat flow
\begin{equation}
 \partial_tF_a=\tau_{\mathbf{g}_a(t),\bar{\mathbf{g}}}(F_a),
 \qquad F_a(0)=\operatorname{id}_M.
 \label{eq:flujo-armonico-aplicaciones-unicidad}
\end{equation}
Corollary~\ref{cor:flujo-armonico-aplicaciones-global} gives a smooth
solution on a common interval. Since
$dF_a(0)=\operatorname{Id}$, uniform continuity of $dF_a(t)$ and
compactness of $M$ allow the interval to be reduced so that every
$dF_a(t)$ is a linear isomorphism. Thus $F_a(t)$ is a local
diffeomorphism.

Let us also verify global injectivity. Choose a finite family
of charts with three concentric coordinate balls
$V_\alpha\Subset W_\alpha\Subset U_\alpha$, such that the
$V_\alpha$ cover $M$ and each $W_\alpha$ is convex in coordinates.
Convergence $F_a(t)\to\operatorname{id}_M$ in $C^1$ allows time to be reduced
so that $F_a(t)(W_\alpha)\subseteq U_\alpha$ and the coordinate
representation $f_{a,t,\alpha}$ satisfies
\[
 \|D(f_{a,t,\alpha}-\operatorname{id})\|_\infty<\frac12
\]
on $W_\alpha$, simultaneously for all indices of the finite
family. Integration along the segment between two points $x,y$
of this ball gives
\[
 |f_{a,t,\alpha}(x)-f_{a,t,\alpha}(y)|
 \geq\frac12|x-y|.
\]
Thus each restriction to $W_\alpha$ is injective.

Let $\delta>0$ be a Lebesgue number of the cover
$(V_\alpha)$ for a fixed Riemannian distance. Reduce time
again so that
$d(F_a(t)(x),x)<\delta/3$ for every $x\in M$. If
$F_a(t)(x)=F_a(t)(y)$, the triangle inequality gives
$d(x,y)<2\delta/3$, so both points belong to some
$V_\alpha$. The local injectivity just proved implies $x=y$.
Finally, $F_a(t)$ preserves each connected component, since
$s\mapsto F_a(st)(x)$ is a path starting at $x$. On each component,
its image is nonempty, open because it is a local diffeomorphism, and closed
by compactness. Connectedness makes it the whole component. Thus $F_a(t)$ is
bijective, and its inverse is smooth by the inverse function theorem.

Define
\[
 \widehat{\mathbf{g}}_a(t)=(F_a(t)^{-1})^*\mathbf{g}_a(t).
\]
The tension field of the identity from $(M,\widehat{\mathbf{g}}_a)$ to
$(M,\bar{\mathbf{g}})$ is
\begin{equation}
 \tau_{\widehat{\mathbf{g}}_a,\bar{\mathbf{g}}}(\operatorname{id})^k
 =\sum_{i,j=1}^m\widehat g_a^{ij}
 \bigl(\overline\Gamma^k_{ij}-\Gamma(\widehat{\mathbf{g}}_a)^k_{ij}\bigr)
 =-\mathbf{W}_{\bar{\mathbf{g}}}(\widehat{\mathbf{g}}_a)^k.
 \label{eq:tension-identidad-deturck}
\end{equation}
Set
\[
 \mathbf{V}_a=(\partial_tF_a)\circ F_a^{-1}.
\]
Precomposing the equation for $F_a$ with $F_a^{-1}$, naturality of the
tension field under changes of domain gives the exact identity
\begin{equation}
 \mathbf{V}_a
 =\tau_{\widehat{\mathbf{g}}_a,\bar{\mathbf{g}}}(\operatorname{id})
 =-\mathbf{W}_{\bar{\mathbf{g}}}(\widehat{\mathbf{g}}_a).
 \label{eq:velocidad-flujo-armonico-deturck}
\end{equation}
Differentiation of pullback by the inverse family gives
\begin{equation}
 \partial_t\widehat{\mathbf{g}}_a
 =(F_a^{-1})^*(\partial_t\mathbf{g}_a)
 -\mathcal L_{\mathbf{V}_a}\widehat{\mathbf{g}}_a.
 \label{eq:derivada-pullback-inverso-flujo-armonico}
\end{equation}
Indeed, differentiating $F_a^{-1}\circ F_a=\operatorname{id}$ identifies
the Eulerian velocity of the inverse with the pullback of $-\mathbf{V}_a$, and rule
\eqref{eq:derivada-pullback-tiempo} gives
\eqref{eq:derivada-pullback-inverso-flujo-armonico}. By
\eqref{eq:naturalidad-ricci} and
\eqref{eq:velocidad-flujo-armonico-deturck}, this formula becomes
\[
 \partial_t\widehat{\mathbf{g}}_a
 =-2\operatorname{Ric}_{\widehat{\mathbf{g}}_a}
  +\mathcal L_{\mathbf{W}_{\bar{\mathbf{g}}}(\widehat{\mathbf{g}}_a)}
   \widehat{\mathbf{g}}_a,
 \qquad \widehat{\mathbf{g}}_a(0)=\mathbf{g}_0.
\]
Uniqueness of the quasilinear parabolic problem gives
$\widehat{\mathbf{g}}_1=\widehat{\mathbf{g}}_2=: \widehat{\mathbf{g}}$. Moreover,
\[
 \partial_tF_a
 =-\mathbf{W}_{\bar{\mathbf{g}}}(\widehat{\mathbf{g}})\circ F_a,
 \qquad F_a(0)=\operatorname{id}_M.
\]
Uniqueness for this nonautonomous ordinary differential equation implies
$F_1=F_2$. Since $\mathbf{g}_a=F_a^*\widehat{\mathbf{g}}$, we conclude that
$\mathbf{g}_1=\mathbf{g}_2$ on
an initial interval. The set of times up to which they agree is
closed by continuity and open by the same argument restarted at
any time; connectedness of the time interval shows that they agree on every common interval of
existence.

For smooth families of smooth data, keep a reference
metric fixed and reduce time uniformly. The fourth assertion of
Corollary~\ref{cor:regularidad-dependencia-parabolica}, parameter dependence of
equation \eqref{eq:flujo-difeomorfismos-deturck}, and the
pullback rule give joint smoothness. Parabolic smoothing allows
it to be measured in every $C^k$ norm on $M\times[\tau,T]$ when $\tau>0$.
\end{proof}

\begin{remark}
The second half of the proof also gives the reverse direction of the
DeTurck procedure. Harmonic map heat flow takes any
Ricci solution to a solution of the same strictly parabolic system;
the ordinary differential equation for the diffeomorphisms recovers the geometric
solution unambiguously.
\end{remark}

\section{Evolution equations}
\label{sec:evolucion-curvatura-ricci}

Let $\mathbf{g}(t)$ be a smooth solution of \eqref{eq:flujo-ricci}. Substituting
$\mathbf{h}=-2\operatorname{Ric}$ into
Lemma~\ref{lem:identidades-variacion-metrica} gives the identities
\begin{align}
 \partial_tg^{ij}&=2\operatorname{Ric}^{ij},
 \label{eq:evolucion-inversa-ricci}\\
 \partial_td\lambda_{\mathbf{g}}&=-R\,d\lambda_{\mathbf{g}},
 \label{eq:evolucion-medida-ricci}\\
 \partial_t\Gamma^k_{ij}
 &=-\sum_{\ell=1}^m g^{k\ell}\bigl(
 \nabla_i\operatorname{Ric}_{j\ell}
 +\nabla_j\operatorname{Ric}_{i\ell}
 -\nabla_\ell\operatorname{Ric}_{ij}\bigr).
 \label{eq:evolucion-conexion-ricci}
\end{align}
In particular,
\begin{equation}
 \frac{d}{dt}\operatorname{Vol}(M,\mathbf{g}(t))
 =-\int_M R_{\mathbf{g}(t)}\,d\lambda_{\mathbf{g}(t)}.
 \label{eq:evolucion-volumen-total-ricci}
\end{equation}

\begin{proposition}[Variation of the Laplacian]
\label{prop:variacion-laplaciano-metrica-tiempo}
Let $M$ be a smooth manifold with or without boundary, let $\mathbf{g}(t)$ be a
smooth family of metrics on $M$, let
$\mathbf{h}=\partial_t\mathbf{g}$, and let $f=f(x,t)$ be a smooth function. Then
\begin{equation}
\partial_t(\Delta f)
=\Delta(\partial_tf)
-\sum_{i,j=1}^m h^{ij}(\nabla^2f)_{ij}
-\sum_{k=1}^m\left(
\sum_{i=1}^m\nabla^ih_i{}^k-\frac12\nabla^k\operatorname{tr}_{\mathbf{g}}\mathbf{h}
\right)\nabla_kf.
\label{eq:variacion-laplaciano-metrica-tiempo}
\end{equation}
Under Ricci flow, the contracted Bianchi identity gives
\begin{equation}
\partial_t(\Delta f)
=\Delta(\partial_tf)
+2\sum_{i,j=1}^m\operatorname{Ric}^{ij}(\nabla^2f)_{ij}.
\label{eq:evolucion-laplaciano-ricci}
\end{equation}
Under a flow of the form
$\partial_t\mathbf{g}=-2\operatorname{Ric}_{\mathbf{g}}+\frac{2}{m}r(t)\mathbf{g}$, with $r(t)$
constant in space, there is the additional term
$-\frac{2}{m}r(t)\,\Delta_{\mathbf{g}} f$.
\end{proposition}

\begin{proof}
In coordinates,
\[
\Delta f=\sum_{i,j=1}^m g^{ij}\left(\partial_i\partial_jf
-\sum_{k=1}^m\Gamma^k_{ij}\partial_kf\right).
\]
Differentiate, use $\partial_tg^{ij}=-h^{ij}$, and contract the Christoffel
variation formula. This gives
\[
\sum_{i,j=1}^m g^{ij}\partial_t\Gamma^k_{ij}
=\sum_{i=1}^m\nabla^ih_i{}^k-\frac12\nabla^k\operatorname{tr}_{\mathbf{g}}\mathbf{h},
\]
and hence \eqref{eq:variacion-laplaciano-metrica-tiempo}. For
$\mathbf{h}=-2\operatorname{Ric}$, the equality
$\displaystyle \sum_{i=1}^m\nabla^i\operatorname{Ric}_{ik}=\frac12\nabla_kR$ cancels the first-order
term. Adding $\frac{2}{m}r\mathbf{g}$, its spatial derivative is zero and
its contraction with $\nabla^2f$ is $\frac{2}{m}r\Delta f$.
\end{proof}

Curvature evolution starts with the uncommuted formula obtained
directly from variation of the connection.

\begin{lemma}[Uncommuted curvature evolution]
\label{lem:evolucion-no-conmutada-riemann}
Let $M$ be a smooth manifold of dimension $m$ with or without boundary and let $\mathbf{g}(t)$ be a
smooth solution of Ricci flow on $M$.
The tensor of type $(1,3)$ satisfies
\begin{align}
 \partial_tR_{ijk}{}^\ell
 =\sum_{q=1}^m g^{\ell q}\bigl(&
 -(\nabla^2\operatorname{Ric})_{ij kq}
 -(\nabla^2\operatorname{Ric})_{ik jq}
 +(\nabla^2\operatorname{Ric})_{iq jk}\\
 &+(\nabla^2\operatorname{Ric})_{ji kq}
 +(\nabla^2\operatorname{Ric})_{jk iq}
 -(\nabla^2\operatorname{Ric})_{jq ik}
 \bigr).
 \label{eq:evolucion-riemann-no-conmutada}
\end{align}
For the fully covariant tensor, add the term coming from the
metric:
\begin{equation}
 \partial_tR_{ijk\ell}
 =-2\sum_{q=1}^m\operatorname{Ric}_{\ell q}R_{ijk}{}^q
  +\sum_{q=1}^m g_{\ell q}\partial_tR_{ijk}{}^q.
 \label{eq:evolucion-riemann-covariante-no-conmutada}
\end{equation}
\end{lemma}

\begin{proof}
Equation \eqref{eq:variacion-riemann-conexion} and
\eqref{eq:evolucion-conexion-ricci} give
\[
 \partial_tR_{ijk}{}^\ell
 =\nabla_i(\partial_t\Gamma^\ell_{jk})
  -\nabla_j(\partial_t\Gamma^\ell_{ik}).
\]
Substitute \eqref{eq:evolucion-conexion-ricci} and use $\nabla \mathbf{g}=0$;
expanding the two parentheses produces exactly the six terms of
\eqref{eq:evolucion-riemann-no-conmutada}. Finally, differentiate
$\displaystyle R_{ijk\ell}=\displaystyle\sum_{q=1}^m g_{\ell q}R_{ijk}{}^q$ and use
$\partial_tg_{\ell q}=-2\operatorname{Ric}_{\ell q}$.
\end{proof}

The commutations needed to recognize the Laplacian follow from
\begin{equation}
 ([\nabla_a,\nabla_b]\mathbf{T})_{i_1\dots i_r}
 =-\sum_{\alpha=1}^r\sum_{p=1}^m
 R_{ab i_\alpha}{}^p
 T_{i_1\dots i_{\alpha-1}p i_{\alpha+1}\dots i_r}
 \label{eq:conmutador-tensor-covariante}
\end{equation}
and the second Bianchi identity
\begin{equation}
 \nabla_aR_{bcde}+\nabla_bR_{cade}+\nabla_cR_{abde}=0.
 \label{eq:segunda-bianchi-ricci-flujo}
\end{equation}
Contracting the latter gives
\begin{equation}
 \sum_{i=1}^m\nabla^i\operatorname{Ric}_{ij}=\frac{1}{2}\nabla_jR.
 \label{eq:bianchi-contraida-ricci-flujo}
\end{equation}

In estimate formulas, use the notation $\mathbf{S}*\mathbf{T}$ for a universal finite
sum of contractions of $\mathbf{S}\otimes \mathbf{T}$ using
$\mathbf{g}$ and $\mathbf{g}^{-1}$.
This is Notation~\ref{notacion estrella};
Lemma~\ref{lema estimacion estrella} controls each contraction. The number
of summands and their coefficients depend only on tensor type and
$m$. For covariant tensors of orders $r,q\in\mathbb N_0$, there exists
$C=C(m,r,q)$ such that
\[
 |\mathbf{S}*\mathbf{T}|_{\mathbf{g}}
 \leq C|\mathbf{S}|_{\mathbf{g}}|\mathbf{T}|_{\mathbf{g}}
\]
for fixed tensor types. This notation does not suppress derivatives: it only collects
algebraic contractions completely determined by applying
\eqref{eq:conmutador-tensor-covariante}.

\begin{proposition}[Evolution of Ricci, scalar, and Riemann curvature]
\label{prop:evolucion-curvaturas-ricci}
Let $M$ be a smooth manifold of dimension $m$ with or without boundary and let $\mathbf{g}(t)$ be a
smooth solution of Ricci flow on $M$.
Define
\begin{equation}
 B_{ijk\ell}
 =\sum_{p,q,r,s=1}^m g^{pr}g^{qs}R_{pijq}R_{rk\ell s}.
 \label{eq:definicion-B-curvatura-ricci}
\end{equation}
Then
\begin{align}
 \partial_t\operatorname{Ric}_{ij}
 &=\Delta\operatorname{Ric}_{ij}
 +2\sum_{k,\ell=1}^mR_{kij\ell}\operatorname{Ric}^{k\ell}
 -2\sum_{k=1}^m\operatorname{Ric}_i{}^k\operatorname{Ric}_{kj},
 \label{eq:evolucion-tensor-ricci}\\
 \partial_tR_{\mathbf{g}(t)}
 &=\Delta_{\mathbf{g}(t)}R_{\mathbf{g}(t)}
   +2|\operatorname{Ric}_{\mathbf{g}(t)}|_{\mathbf{g}(t)}^2,
 \label{eq:evolucion-curvatura-escalar}
\end{align}
Moreover, with all indices covariant,
\begin{align}
 \partial_tR_{ijk\ell}
 ={}&\Delta R_{ijk\ell}
 +2\bigl(B_{ij\ell k}-B_{ijk\ell}
 -B_{ikj\ell}+B_{i\ell jk}\bigr)
 \nonumber\\
 &-\sum_{p,q=1}^m g^{pq}\bigl(
 R_{pjk\ell}\operatorname{Ric}_{qi}
 +R_{ipk\ell}\operatorname{Ric}_{qj}
 +R_{ijkp}\operatorname{Ric}_{q\ell}
 +R_{ijp\ell}\operatorname{Ric}_{qk}\bigr).
 \label{eq:evolucion-riemann-expandida}
\end{align}
In particular,
\begin{equation}
 (\partial_t-\Delta)\operatorname{Rm}
 =\operatorname{Rm}*\operatorname{Rm}.
 \label{eq:evolucion-riemann-esquematica}
\end{equation}
The schematic equality is exactly
\eqref{eq:evolucion-riemann-expandida}: the right-hand side collects only
algebraic contractions of two curvature tensors and contains no
additional derivatives.
\end{proposition}

\begin{proof}
First derive the four-index formula. Commuting in
\eqref{eq:evolucion-riemann-no-conmutada} the two derivatives falling on
$\operatorname{Ric}_{kq}$ and using
\eqref{eq:evolucion-riemann-covariante-no-conmutada} gives
\begin{align}
 \partial_tR_{ijk\ell}
 ={}&(\nabla^2\operatorname{Ric})_{i \ell jk}
  -(\nabla^2\operatorname{Ric})_{ik j\ell}
  -(\nabla^2\operatorname{Ric})_{j \ell ik}
  +(\nabla^2\operatorname{Ric})_{jk i\ell}
 \nonumber\\
 &-\sum_{p,q=1}^m g^{pq}\bigl(
 R_{ijkp}\operatorname{Ric}_{q\ell}
 +R_{ijp\ell}\operatorname{Ric}_{qk}\bigr).
 \label{eq:evolucion-riemann-cuatro-indices-intermedia}
\end{align}
Indeed, the four terms written with two derivatives are the four
uncommuted terms of
\eqref{eq:evolucion-riemann-no-conmutada}. The two remaining terms form
\[
 ([\nabla_j,\nabla_i]\operatorname{Ric})_{k\ell}
 =\sum_{p=1}^mR_{ijk}{}^p\operatorname{Ric}_{p\ell}
  +\sum_{p=1}^mR_{ij\ell}{}^p\operatorname{Ric}_{kp},
\]
by \eqref{eq:conmutador-tensor-covariante}. Adding the term coming
from lowering the index gives
\begin{align*}
 &\sum_{p=1}^mR_{ijk}{}^p\operatorname{Ric}_{p\ell}
 +\sum_{p=1}^mR_{ij\ell}{}^p\operatorname{Ric}_{kp}
 -2\sum_{q=1}^m\operatorname{Ric}_{\ell q}R_{ijk}{}^q\\
 &\qquad=-\sum_{p,q=1}^m g^{pq}\bigl(
 R_{ijkp}\operatorname{Ric}_{q\ell}
 +R_{ijp\ell}\operatorname{Ric}_{qk}\bigr),
\end{align*}
where $R_{ij\ell p}=-R_{ijp\ell}$ was used. This proves
\eqref{eq:evolucion-riemann-cuatro-indices-intermedia}.

We now prove the Simons identity transforming the second
derivatives. The second Bianchi identity and its contraction give, with the
index order used here,
\begin{align}
 \Delta R_{ijk\ell}
 &=\sum_{p,q=1}^m g^{pq}(\nabla^2\operatorname{Rm})_{p i qjk\ell}
   -\sum_{p,q=1}^m g^{pq}(\nabla^2\operatorname{Rm})_{p j qik\ell},
 \label{eq:laplaciano-riemann-bianchi}\\
 \sum_{p,q=1}^m g^{pq}\nabla_pR_{qjk\ell}
 &=\nabla_\ell\operatorname{Ric}_{jk}
   -\nabla_k\operatorname{Ric}_{j\ell}.
 \label{eq:bianchi-contraida-cuatro-indices}
\end{align}
To carry out the commutation, fix a point and an orthonormal frame
there. Identity \eqref{eq:conmutador-tensor-covariante} gives
\begin{align*}
 &\sum_{p=1}^m
 ([\nabla_p,\nabla_i]\operatorname{Rm})_{pjk\ell}\\
 &\quad=\sum_{r=1}^m\operatorname{Ric}_{ir}R_{rjk\ell}
 -\sum_{p,r=1}^m\bigl(
 R_{pijr}R_{prk\ell}+R_{pikr}R_{pjr\ell}
                         +R_{pi\ell r}R_{pjkr}\bigr).
\end{align*}
The first term comes from
$\displaystyle \sum_{p=1}^mR_{pipr}=-\operatorname{Ric}_{ir}$.
The first Bianchi identity and antisymmetry in the last slots
transform the three remaining sums into
\begin{align*}
 -\sum_{p,r=1}^mR_{pijr}R_{prk\ell}
 &=B_{ijk\ell}-B_{ij\ell k},\\
 -\sum_{p,r=1}^mR_{pikr}R_{pjr\ell}&=B_{ikj\ell},\\
 -\sum_{p,r=1}^mR_{pi\ell r}R_{pjkr}&=-B_{i\ell jk}.
\end{align*}
For example, the first equality follows by substituting
$R_{prk\ell}=-R_{pk\ell r}+R_{p\ell kr}$, which is Bianchi with
indices permuted. The other two equalities use directly
$R_{pjr\ell}=-R_{pj\ell r}$ and the definition of $B$.
Returning to arbitrary coordinates gives
\begin{align*}
 &\sum_{p,q=1}^m g^{pq}\bigl(
 (\nabla^2\operatorname{Rm})_{p i qjk\ell}
 -(\nabla^2\operatorname{Rm})_{i p qjk\ell}\bigr)\\
 &\qquad=\sum_{r=1}^m\operatorname{Ric}_i{}^rR_{rjk\ell}
 -B_{ij\ell k}+B_{ijk\ell}-B_{i\ell jk}+B_{ikj\ell}.
\end{align*}
Substitute this equality and its analogue with $i$ and $j$ interchanged into
\eqref{eq:laplaciano-riemann-bianchi}, then use
\eqref{eq:bianchi-contraida-cuatro-indices}. Grouping the four terms
$B$ yields
\begin{align}
 &\Delta R_{ijk\ell}
 +2\bigl(B_{ij\ell k}-B_{ijk\ell}
 -B_{ikj\ell}+B_{i\ell jk}\bigr)
 \nonumber\\
 &\quad=(\nabla^2\operatorname{Ric})_{i \ell jk}
 -(\nabla^2\operatorname{Ric})_{ik j\ell}
 -(\nabla^2\operatorname{Ric})_{j \ell ik}
 +(\nabla^2\operatorname{Ric})_{jk i\ell}
 \nonumber\\
 &\qquad+\sum_{p,q=1}^m g^{pq}\bigl(
 R_{pjk\ell}\operatorname{Ric}_{qi}
 +R_{ipk\ell}\operatorname{Ric}_{qj}\bigr).
 \label{eq:identidad-simons-riemann-ricci}
\end{align}
Substituting \eqref{eq:identidad-simons-riemann-ricci} into
\eqref{eq:evolucion-riemann-cuatro-indices-intermedia} proves
\eqref{eq:evolucion-riemann-expandida} term by term.

To obtain the Ricci equation, substitute
$\mathbf{h}=-2\operatorname{Ric}$ into \eqref{eq:linealizacion-ricci}. The contracted
Bianchi identity gives
\[
 \delta_{\mathbf{g}}(-2\operatorname{Ric})=dR,
 \qquad
 \operatorname{tr}_{\mathbf{g}}(-2\operatorname{Ric})=-2R.
\]
Thus the two Hessian terms cancel. The explicit formula
\eqref{eq:termino-cero-linealizacion-ricci} leaves
\[
 \partial_t\operatorname{Ric}_{ij}
 =\Delta\operatorname{Ric}_{ij}
  +2\sum_{k,\ell=1}^mR_{kij\ell}\operatorname{Ric}^{k\ell}
  -2\sum_{k=1}^m\operatorname{Ric}_i{}^k\operatorname{Ric}_{kj},
\]
which is \eqref{eq:evolucion-tensor-ricci}.

The scalar equation can also be obtained without contracting the preceding one.
Indeed, \eqref{eq:variacion-curvatura-escalar}, with
$\mathbf{h}=-2\operatorname{Ric}$, gives
\begin{align*}
 \partial_tR_{\mathbf{g}(t)}
 &=2|\operatorname{Ric}_{\mathbf{g}(t)}|_{\mathbf{g}(t)}^2
 -2\operatorname{div}_{\mathbf{g}(t)}^2\operatorname{Ric}_{\mathbf{g}(t)}
 +2\Delta_{\mathbf{g}(t)}R_{\mathbf{g}(t)}\\
 &=2|\operatorname{Ric}_{\mathbf{g}(t)}|_{\mathbf{g}(t)}^2
   +\Delta_{\mathbf{g}(t)}R_{\mathbf{g}(t)},
\end{align*}
where \eqref{eq:bianchi-contraida-ricci-flujo} was used. This proves
\eqref{eq:evolucion-curvatura-escalar}.

Finally, every additional term in
\eqref{eq:evolucion-riemann-expandida} is a contraction of two copies of
$\operatorname{Rm}$, because $\operatorname{Ric}$ is a contraction of
$\operatorname{Rm}$. This proves
\eqref{eq:evolucion-riemann-esquematica}.
\end{proof}

\begin{corollary}[Evolution of curvature norms]
\label{cor:evolucion-normas-curvatura-ricci}
Let $M$ be a smooth manifold of dimension $m$ with or without boundary and let $\mathbf{g}(t)$ be a
smooth solution of Ricci flow on $M$.
There exist constants $C_m,C'_m>0$, depending only on dimension, such
that
\begin{align}
 (\partial_t-\Delta_{\mathbf{g}(t)})
 |\operatorname{Rm}_{\mathbf{g}(t)}|_{\mathbf{g}(t)}^2
 &=-2|\nabla_{\mathbf{g}(t)}\operatorname{Rm}_{\mathbf{g}(t)}|_{\mathbf{g}(t)}^2
   -8\sum_{i,j,k,\ell=1}^m\bigl(B_{ijk\ell}+B_{ikj\ell}\bigr)R^{ijk\ell},
 \label{eq:evolucion-norma-riemann}\\
 (\partial_t-\Delta_{\mathbf{g}(t)})
 |\operatorname{Rm}_{\mathbf{g}(t)}|_{\mathbf{g}(t)}^2
 &\leq -2|\nabla_{\mathbf{g}(t)}\operatorname{Rm}_{\mathbf{g}(t)}|_{\mathbf{g}(t)}^2
   +C_m|\operatorname{Rm}_{\mathbf{g}(t)}|_{\mathbf{g}(t)}^3,
 \label{eq:desigualdad-evolucion-norma-riemann}\\
 (\partial_t-\Delta_{\mathbf{g}(t)})
 |\operatorname{Ric}_{\mathbf{g}(t)}|_{\mathbf{g}(t)}^2
 &=-2|\nabla_{\mathbf{g}(t)}\operatorname{Ric}_{\mathbf{g}(t)}|_{\mathbf{g}(t)}^2
   +4\sum_{i,j,k,\ell=1}^m R_{kij\ell}\operatorname{Ric}^{ij}\operatorname{Ric}^{k\ell},
 \label{eq:evolucion-norma-ricci}\\
 (\partial_t-\Delta_{\mathbf{g}(t)})R_{\mathbf{g}(t)}^2
 &=-2|dR_{\mathbf{g}(t)}|_{\mathbf{g}(t)}^2
   +4R_{\mathbf{g}(t)}|\operatorname{Ric}_{\mathbf{g}(t)}|_{\mathbf{g}(t)}^2.
 \label{eq:evolucion-cuadrado-escalar}
\end{align}
In particular,
\begin{equation}
 (\partial_t-\Delta_{\mathbf{g}(t)})
 |\operatorname{Ric}_{\mathbf{g}(t)}|_{\mathbf{g}(t)}^2
 \leq -2|\nabla_{\mathbf{g}(t)}\operatorname{Ric}_{\mathbf{g}(t)}|_{\mathbf{g}(t)}^2
 +C'_m|\operatorname{Rm}_{\mathbf{g}(t)}|_{\mathbf{g}(t)}
 |\operatorname{Ric}_{\mathbf{g}(t)}|_{\mathbf{g}(t)}^2.
 \label{eq:desigualdad-evolucion-norma-ricci}
\end{equation}
\end{corollary}

\begin{proof}
For any time-dependent covariant tensor $\mathbf{T}$,
\[
 \Delta_{\mathbf{g}(t)}|\mathbf{T}|_{\mathbf{g}(t)}^2
 =2\langle\Delta_{\mathbf{g}(t)}\mathbf{T},\mathbf{T}\rangle_{\mathbf{g}(t)}
  +2|\nabla_{\mathbf{g}(t)}\mathbf{T}|_{\mathbf{g}(t)}^2.
\]
Time differentiation of the norm also differentiates each copy of $\mathbf{g}^{-1}$.
By \eqref{eq:evolucion-inversa-ricci}, these terms are contractions of
$\operatorname{Ric}*\mathbf{T}*\mathbf{T}$. With
$\mathbf{T}=\operatorname{Rm}$, contract
\eqref{eq:evolucion-riemann-expandida} with $2R^{ijk\ell}$ and add the
four terms obtained by differentiating the inverse metrics. The terms
$\operatorname{Ric}*\operatorname{Rm}*\operatorname{Rm}$ cancel and,
for every covariant slot, the Ricci term in the tensor evolution
cancels the term differentiating the inverse metric in that same
slot. For the summands $B$, renaming $k$ and $\ell$ and using
$R^{ij\ell k}=-R^{ijk\ell}$ gives
\begin{align*}
 \sum_{i,j,k,\ell=1}^m B_{ij\ell k}R^{ijk\ell}
 &=-\sum_{i,j,k,\ell=1}^m B_{ijk\ell}R^{ijk\ell},\\
 \sum_{i,j,k,\ell=1}^m B_{i\ell jk}R^{ijk\ell}
 &=-\sum_{i,j,k,\ell=1}^m B_{ikj\ell}R^{ijk\ell}.
\end{align*}
Consequently, the remainder is
\[
-8\sum_{i,j,k,\ell=1}^m(B_{ijk\ell}+B_{ikj\ell})R^{ijk\ell}.
\]
This proves \eqref{eq:evolucion-norma-riemann}. The algebraic estimate for
$B$ gives
\eqref{eq:desigualdad-evolucion-norma-riemann}.

For the Ricci tensor, directly differentiate the two inverse metrics
occurring in
$\displaystyle |\operatorname{Ric}|^2=\displaystyle\sum_{i,j,a,b=1}^m g^{ia}g^{jb}\operatorname{Ric}_{ij}
\operatorname{Ric}_{ab}$ and use
\eqref{eq:evolucion-inversa-ricci}. This yields
\[
 \partial_t|\operatorname{Ric}_{\mathbf{g}(t)}|_{\mathbf{g}(t)}^2
 =2\langle\partial_t\operatorname{Ric}_{\mathbf{g}(t)},
 \operatorname{Ric}_{\mathbf{g}(t)}\rangle_{\mathbf{g}(t)}
  +4\sum_{i,j,k=1}^m\operatorname{Ric}_i{}^j\operatorname{Ric}_j{}^k
    \operatorname{Ric}_k{}^i.
\]
The cubic Ricci terms coming from
\eqref{eq:evolucion-tensor-ricci} cancel the second summand exactly.
Subtracting the formula for
$\Delta_{\mathbf{g}(t)}|\operatorname{Ric}_{\mathbf{g}(t)}|_{\mathbf{g}(t)}^2$ leaves
\eqref{eq:evolucion-norma-ricci}. The following inequality is Cauchy--Schwarz
in each fiber. Finally,
\[
 (\partial_t-\Delta_{\mathbf{g}(t)})R_{\mathbf{g}(t)}^2
 =2R_{\mathbf{g}(t)}(\partial_t-\Delta_{\mathbf{g}(t)})R_{\mathbf{g}(t)}
  -2|dR_{\mathbf{g}(t)}|_{\mathbf{g}(t)}^2,
\]
and \eqref{eq:evolucion-curvatura-escalar} proves
\eqref{eq:evolucion-cuadrado-escalar}.
\end{proof}

\section{The maximum principle and a priori estimates}
\label{sec:maximo-estimaciones-flujo-ricci}

The evolution equations become useful when combined with the
maximum principle. Scalar curvature satisfies a closed inequality,
while the norm of the Riemann tensor obeys an equation with a cubic
term. These two observations give the first bounds and, through
a Bernstein induction, estimates for all derivatives of
curvature.

To apply the maximum principle from the preceding chapter, fix
$\nabla^0=\nabla^{\mathbf{g}(0)}$. At every time, on functions,
\[
 \Delta_{\mathbf{g}(t)}f
 =\sum_{i,j=1}^m g^{ij}(t)((\nabla^0)^2f)_{ij}
    +\sum_{k=1}^m b^k(x,t)\nabla_k^0 f,
\]
where
\[
 b^k(x,t)=\sum_{i,j=1}^m g^{ij}(t)
 \bigl(\Gamma(\mathbf{g}(0))^k_{ij}-\Gamma(\mathbf{g}(t))^k_{ij}\bigr),
 \qquad k\in\{1,\ldots,m\}.
\]
On a compact subinterval of existence, smoothness and positivity of
$\mathbf{g}(t)$ make ellipticity and the coefficient bounds
uniform. Thus the operator satisfies the hypotheses of
Theorem~\ref{teo:principio-maximo-parabolico-escalar}, even though the metric
defining the Laplacian varies with time.

\begin{proposition}[Lower bound for scalar curvature]
\label{prop:cota-inferior-escalar-ricci}
Let $M$ be a closed smooth manifold of dimension $m$ and let
$\mathbf{g}(t)$ be a smooth solution of Ricci flow on $M$.
Let
\[
 a_0=\min_{x\in M}R_{\mathbf{g}(0)}(x).
\]
As long as the denominator is positive,
\begin{equation}
 R_{\mathbf{g}(t)}(x)\geq
 \frac{a_0}{1-\frac{2}{m}a_0t}
 \qquad (x\in M).
 \label{eq:cota-inferior-escalar-ricci}
\end{equation}
In particular, the condition $R_{\mathbf{g}(0)}\geq0$ is preserved. If
$R_{\mathbf{g}(0)}>0$, then $R_{\mathbf{g}(t)}>0$ throughout the existence interval.
\end{proposition}

\begin{proof}
Cauchy--Schwarz for the Ricci eigenvalues gives
\[
 |\operatorname{Ric}_{\mathbf{g}(t)}|_{\mathbf{g}(t)}^2
 \geq\frac{1}{m}R_{\mathbf{g}(t)}^2.
\]
Let $a(t)$ be the solution of
\[
 a'(t)=\frac{2}{m}a(t)^2,
 \qquad a(0)=a_0.
\]
If $a_0\neq0$, then
\[
 a(t)=\frac{a_0}{1-\frac{2}{m}a_0t};
\]
if $a_0=0$, we have $a(t)=0$. For $u=R_{\mathbf{g}(t)}-a$,
\begin{align*}
 (\partial_t-\Delta)u
 &=2|\operatorname{Ric}_{\mathbf{g}(t)}|_{\mathbf{g}(t)}^2-\frac{2}{m}a^2\\
 &\geq\frac{2}{m}(R_{\mathbf{g}(t)}^2-a^2)
 =\frac{2}{m}(R_{\mathbf{g}(t)}+a)u.
\end{align*}
Set $b=\frac{2}{m}(R_{\mathbf{g}(t)}+a)$ and $v=-u$. On every compact subinterval
of existence, $b$ has an upper bound $B$. The preceding inequality
reads
\[
 (\partial_t-\Delta)v\leq bv.
\]
For $\widetilde v=e^{-Bt}v$ we have
\[
 (\partial_t-\Delta)\widetilde v
 \leq (b-B)\widetilde v.
\]
Since $b-B\leq0$ and $\widetilde v(0)\leq0$,
Theorem~\ref{teo:principio-maximo-parabolico-escalar} applied to the positive
maximum of $\widetilde v$ gives $\widetilde v\leq0$. Thus $u\geq0$.
Strict positivity also follows from the principle already proved:
if $\displaystyle a_0=\min_{x\in M}R_{\mathbf{g}(0)}(x)>0$, then
$(\partial_t-\Delta)(a_0-R)=-2|\operatorname{Ric}|^2\leq0$.
The maximum principle implies $R(x,t)\geq a_0>0$ on every compact
subinterval of existence. The same argument with $a_0=0$ preserves
nonnegativity.
\end{proof}

\begin{proposition}[Initial curvature control]
\label{prop:control-inicial-curvatura-ricci}
Let $M$ be a closed smooth manifold of dimension $m$ and let
$\mathbf{g}(t)$ be a smooth solution of Ricci flow on $M$.
If
\[
 K_0=\max_{x\in M}|\operatorname{Rm}_{\mathbf{g}(0)}(x)|_{\mathbf{g}(0)},
\]
then there exists $c_m>0$ such that
\begin{equation}
 \max_{x\in M}|\operatorname{Rm}_{\mathbf{g}(t)}(x)|_{\mathbf{g}(t)}
 \leq\frac{K_0}{1-c_mK_0t}
 \label{eq:control-inicial-curvatura-ricci}
\end{equation}
if $1-c_mK_0t>0$. When $K_0=0$, the flow remains flat.
\end{proposition}

\begin{proof}
Set $q(x,t)=|\operatorname{Rm}_{\mathbf{g}(t)}(x)|_{\mathbf{g}(t)}^2$.
Equation~\eqref{eq:desigualdad-evolucion-norma-riemann} and the scalar maximum
principle show that the upper right derivative of the maximum
$\displaystyle Q(t)=\displaystyle\max_{x\in M}q(x,t)$ satisfies
\[
 D^+Q(t)\leq C_mQ(t)^{\frac{3}{2}}.
\]
The solution of $y'=C_my^{\frac{3}{2}}$, $y(0)=K_0^2$, is
\[
 y(t)=\frac{K_0^2}{\left(1-\frac{C_m}{2}K_0t\right)^2}.
\]
Comparison for ordinary differential equations gives the conclusion with
$c_m=\frac{C_m}{2}$. If $K_0=0$, comparison gives $q=0$.
\end{proof}

Curvature derivative estimates are needed to
rule out loss of regularity while curvature remains bounded.

\begin{proposition}[Global derivative estimates]
\label{prop:estimaciones-derivadas-curvatura-ricci}
Let $M$ be a closed smooth manifold of dimension $m$ and let
$\mathbf{g}(t)$ be a smooth solution of Ricci flow on $M$.
Suppose that $\mathbf{g}(t)$ is defined on $M\times[0,T]$ and
\[
 |\operatorname{Rm}_{\mathbf{g}(t)}(x)|_{\mathbf{g}(t)}\leq K
 \qquad ((x,t)\in M\times[0,T]).
\]
For every integer $r\geq1$ there exists a constant
$C_r=C_r(m,r,KT)$ such that
\begin{equation}
 |\nabla_{\mathbf{g}(t)}^r\operatorname{Rm}_{\mathbf{g}(t)}(x)|_{\mathbf{g}(t)}
 \leq C_rK t^{-\frac{r}{2}}
 \qquad ((x,t)\in M\times(0,T]).
 \label{eq:estimaciones-derivadas-curvatura-ricci}
\end{equation}
If $K=0$, all derivatives vanish and the assertion is interpreted
accordingly.
\end{proposition}

\begin{proof}
Write $\mathcal D=\partial_t-\Delta$ and
$\mathbf{T}_r=\nabla^r\operatorname{Rm}$ for $r\in\mathbb N_0$.
If $\mathbf{T}$ is covariant of order $q$, differentiating its components gives
\[
 ([\partial_t,\nabla]\mathbf{T})_{i j_1\ldots j_q}
 =-\sum_{\alpha=1}^q\sum_{p=1}^m
 (\partial_t\Gamma^p_{i j_\alpha})
 T_{j_1\ldots j_{\alpha-1}p j_{\alpha+1}\ldots j_q}.
\]
By \eqref{eq:evolucion-conexion-ricci}, this is a sum of
contractions of $\nabla\operatorname{Rm}$ with $\mathbf{T}$.
On the other hand, Lemma~\ref{lema conmutador laplaciano orden m}, applied
to the tensor connection with $\Delta=-\nabla^*\nabla$, gives
\[
 [\Delta,\nabla]\mathbf{T}
 =\operatorname{Rm}*\nabla\mathbf{T}
       +\nabla\operatorname{Rm}*\mathbf{T}.
\]
Here the bundle in the lemma is the bundle of tensors of the type of $\mathbf{T}$;
its curvature is the sum of the Riemann actions on each slot, as
expressed by \eqref{eq:conmutador-tensor-covariante}.

For every $r\in\mathbb N_0$, we prove the identity
\begin{equation}
 (\partial_t-\Delta_{\mathbf{g}(t)})
 \nabla_{\mathbf{g}(t)}^r\operatorname{Rm}_{\mathbf{g}(t)}
 =\sum_{\substack{a,b\in\mathbb N_0\\a+b=r}}
 \nabla_{\mathbf{g}(t)}^a\operatorname{Rm}_{\mathbf{g}(t)}
 *\nabla_{\mathbf{g}(t)}^b\operatorname{Rm}_{\mathbf{g}(t)}.
 \label{eq:evolucion-derivadas-riemann}
\end{equation}
For $r=0$, this is \eqref{eq:evolucion-riemann-esquematica}. Suppose that
it holds for some $r\geq0$. Since $\mathbf{T}_{r+1}=\nabla\mathbf{T}_r$,
\[
 \mathcal D\mathbf{T}_{r+1}
 =\nabla(\mathcal D\mathbf{T}_r)
       +[\partial_t,\nabla]\mathbf{T}_r-[\Delta,\nabla]\mathbf{T}_r.
\]
Proposition~\ref{prop: leibniz operador *} turns the first term
into the sum, for $a+b=r$, of
$\mathbf{T}_{a+1}*\mathbf{T}_b$ and $\mathbf{T}_a*\mathbf{T}_{b+1}$.
The two commutators contribute only
$\mathbf{T}_1*\mathbf{T}_r$ and $\mathbf{T}_0*\mathbf{T}_{r+1}$.
All these products have total order $r+1$; grouping those
sharing the same pair of orders gives
\eqref{eq:evolucion-derivadas-riemann} for $r+1$. The number of terms
and their coefficients depend only on $m$ and $r$.

Passing to the norm of $\mathbf{T}_r$ differentiates each of the $r+4$
copies of the inverse metric. Identity
$\partial_t\mathbf{g}^{-1}=2\operatorname{Ric}^{\sharp\sharp}$ bounds their
contribution by $C_r|\operatorname{Rm}|\,|\mathbf{T}_r|^2$, which is the
type of term corresponding to $a=0$ or $b=0$ in the following sum.
By the Leibniz rule,
\begin{align}
 (\partial_t-\Delta_{\mathbf{g}(t)})
 |\nabla_{\mathbf{g}(t)}^r\operatorname{Rm}_{\mathbf{g}(t)}|_{\mathbf{g}(t)}^2
 \leq{}&-2|\nabla_{\mathbf{g}(t)}^{r+1}\operatorname{Rm}_{\mathbf{g}(t)}|_{\mathbf{g}(t)}^2\\
 &+C_r\sum_{\substack{a,b\in\mathbb N_0\\a+b=r}}
 |\nabla_{\mathbf{g}(t)}^r\operatorname{Rm}_{\mathbf{g}(t)}|_{\mathbf{g}(t)}
 |\nabla_{\mathbf{g}(t)}^a\operatorname{Rm}_{\mathbf{g}(t)}|_{\mathbf{g}(t)}
 |\nabla_{\mathbf{g}(t)}^b\operatorname{Rm}_{\mathbf{g}(t)}|_{\mathbf{g}(t)}.
 \label{eq:desigualdad-derivadas-riemann}
\end{align}

Scale normalization avoids mixing different powers of $K$. If
$K>0$, define
\begin{equation}
 \widetilde{\mathbf{g}}(\tau)=K \mathbf{g}\left(\frac{\tau}{K}\right),
 \qquad 0\leq\tau\leq KT.
 \label{eq:reescalamiento-prueba-shi}
\end{equation}
Then $\widetilde{\mathbf{g}}$ satisfies Ricci flow and
$|\operatorname{Rm}_{\widetilde{\mathbf{g}}(\tau)}|_{\widetilde{\mathbf{g}}(\tau)}\leq1$.
In the original
variables, the term with no derivatives in
\eqref{eq:desigualdad-evolucion-norma-riemann} has size $K^3$, not
$K^2$; transformation
\eqref{eq:reescalamiento-prueba-shi} turns it into a dimensionless
constant. It therefore suffices to prove the estimate when $K=1$ on an
interval $[0,S]$, with $S=KT$.

Write
\[
 u_j(x,\tau)
 =|\nabla_{\widetilde{\mathbf{g}}(\tau)}^j
 \operatorname{Rm}_{\widetilde{\mathbf{g}}(\tau)}(x)|_{\widetilde{\mathbf{g}}(\tau)}^2.
\]
Proceed by induction on $r$. Suppose the bounds
$u_j\leq C_{j,S}\tau^{-j}$ are known for $j<r$; for $r=1$ this hypothesis uses only
$u_0\leq1$. In \eqref{eq:desigualdad-derivadas-riemann}, summands with
$a=0$ or $b=0$ are bounded by $C_ru_r$. If $1\leq a,b\leq r-1$, the
induction hypothesis gives
$\tau^r\sqrt{u_ru_au_b}\leq C_{r,S}\tau^{r/2}\sqrt{u_r}$.
Write the last product as
$C_{r,S}(\tau^{(r-1)/2}\sqrt{u_r})\tau^{1/2}$ and apply Young.
Since $\tau\leq S$, we obtain
\[
 \tau^r u_r^{\frac{1}{2}}u_a^{\frac{1}{2}}u_b^{\frac{1}{2}}
 \leq \varepsilon\tau^{r-1}u_r+C_{r,S,\varepsilon},
 \qquad a+b=r.
\]
Differentiating the power of $\tau$ and also absorbing
$\tau^ru_r\leq S\tau^{r-1}u_r$ gives
\begin{equation}
 (\partial_\tau-\Delta_{\widetilde{\mathbf{g}}(\tau)})(\tau^ru_r)
 \leq-2\tau^ru_{r+1}+C_{r,S}\tau^{r-1}u_r+C_{r,S}.
 \label{eq:shi-peso-orden-r}
\end{equation}
The same calculation, already justified by the induction hypothesis, gives for
$1\leq j<r$
\begin{equation}
 (\partial_\tau-\Delta_{\widetilde{\mathbf{g}}(\tau)})(\tau^ju_j)
 \leq-2\tau^ju_{j+1}+C_{j,S}\tau^{j-1}u_j+C_{j,S},
 \label{eq:shi-pesos-inferiores}
\end{equation}
and for $j=0$ gives
\begin{equation}
 (\partial_\tau-\Delta_{\widetilde{\mathbf{g}}(\tau)})u_0
 \leq-2u_1+C_m.
 \label{eq:shi-peso-cero}
\end{equation}

Set $A_r=1$ and, for $j=r,r-1,\ldots,1$, define
$A_{j-1}=1+C_{j,S}A_j$. This choice makes the term
$C_{j,S}A_j\tau^{j-1}u_j$ in
\eqref{eq:shi-peso-orden-r} or
\eqref{eq:shi-pesos-inferiores} absorbable by
$-2A_{j-1}\tau^{j-1}u_j$. For
\begin{equation}
 F_r=\tau^ru_r+\sum_{j=0}^{r-1}A_j\tau^ju_j
 \label{eq:funcional-inductivo-shi}
\end{equation}
we therefore obtain
\begin{equation}
 (\partial_\tau-\Delta_{\widetilde{\mathbf{g}}(\tau)})F_r\leq C_{r,S}.
 \label{eq:desigualdad-funcional-inductivo-shi}
\end{equation}
Since $F_r(0)=A_0u_0(0)\leq A_0$, the maximum principle implies
$F_r\leq A_0+C_{r,S}S$. In particular,
\[
 |\nabla_{\widetilde{\mathbf{g}}(\tau)}^r
 \operatorname{Rm}_{\widetilde{\mathbf{g}}(\tau)}|_{\widetilde{\mathbf{g}}(\tau)}
 \leq C_{r,S}\tau^{-\frac{r}{2}}.
\]
This closes the induction.

Finally, under the constant rescaling in
\eqref{eq:reescalamiento-prueba-shi},
\[
 |\nabla_{\widetilde{\mathbf{g}}(\tau)}^r
 \operatorname{Rm}_{\widetilde{\mathbf{g}}(\tau)}|_{\widetilde{\mathbf{g}}(\tau)}
 =K^{-1-\frac{r}{2}}
  |\nabla_{\mathbf{g}(\tau/K)}^r
  \operatorname{Rm}_{\mathbf{g}(\tau/K)}|_{\mathbf{g}(\tau/K)}.
\]
Taking $\tau=Kt$ gives
\[
 |\nabla_{\mathbf{g}(t)}^r\operatorname{Rm}_{\mathbf{g}(t)}|_{\mathbf{g}(t)}
 \leq C_r(m,r,KT)K t^{-\frac{r}{2}},
\]
which is \eqref{eq:estimaciones-derivadas-curvatura-ricci}. The case $K=0$
already follows from Proposition~\ref{prop:control-inicial-curvatura-ricci}.
\end{proof}

\begin{theorem}[Continuation criterion]
\label{teo:criterio-continuacion-flujo-ricci}
Let $\mathbf{g}(t)$ be the maximal solution of Ricci flow on a closed smooth manifold,
defined on $[0,T_{\displaystyle\max})$. If $T_{\displaystyle\max}<\infty$, then
\begin{equation}
 \lim_{t\to T_{\max}^{-}}
 \max_{x\in M}|\operatorname{Rm}_{\mathbf{g}(t)}(x)|_{\mathbf{g}(t)}=\infty.
 \label{eq:explosion-curvatura-tiempo-maximal}
\end{equation}
\end{theorem}

\begin{proof}
Suppose, toward a contradiction, that
$|\operatorname{Rm}_{\mathbf{g}(t)}|_{\mathbf{g}(t)}\leq K$ on
$M\times[0,T_{\displaystyle\max})$. Since
$|\operatorname{Ric}_{\mathbf{g}(t)}|_{\mathbf{g}(t)}
\leq\sqrt m|\operatorname{Rm}_{\mathbf{g}(t)}|_{\mathbf{g}(t)}$, for every vector
$v\neq0$ fixed in $TM$ we have
\[
 \left|\partial_t\log \mathbf{g}(t)(v,v)\right|
 =2\frac{|\operatorname{Ric}_{\mathbf{g}(t)}(v,v)|}{\mathbf{g}(t)(v,v)}
 \leq2\sqrt m K.
\]
Integration in time gives
\begin{equation}
 e^{-2\sqrt mKT_{\max}}\mathbf{g}(0)
 \leq \mathbf{g}(t)\leq
 e^{2\sqrt mKT_{\max}}\mathbf{g}(0).
 \label{eq:equivalencia-uniforme-metricas-ricci}
\end{equation}
Inequality \eqref{eq:equivalencia-uniforme-metricas-ricci}
states that, for every $x\in M$ and every $v\in T_xM$,
\[
 e^{-2\sqrt mKT_{\max}}\mathbf g(0)_x(v,v)
 \leq\mathbf g(t)_x(v,v)
 \leq e^{2\sqrt mKT_{\max}}\mathbf g(0)_x(v,v).
\]
The constants are independent of $x$ and $t<T_{\displaystyle\max}$.
In particular, any possible limit cannot lose positivity.
Proposition~\ref{prop:comparacion-metricas-duales-tensores-volumen}
allows every tensor of fixed type to be measured with $\mathbf g(0)$
or $\mathbf g(t)$ using common constants. To replace
also the $t$-dependent Levi--Civita derivatives by
derivatives with respect to a fixed connection, we need the following
change-of-connection argument.

Apply Proposition~\ref{prop:estimaciones-derivadas-curvatura-ricci}
to the shifted flow starting at $\frac{T_{\displaystyle\max}}{4}$. All covariant derivatives
of Riemann are uniformly bounded on
$M\times[\frac{T_{\displaystyle\max}}{2},T_{\displaystyle\max})$.

Let us explain the passage to a fixed connection. Let $\nabla^0=\nabla^{\mathbf{g}(0)}$ and
\[
 \mathbf{A}(t)=\nabla^{\mathbf{g}(t)}-\nabla^0.
\]
The change-of-connection identities have the form
\begin{equation}
 \nabla^0\mathbf{T}=\nabla^{\mathbf{g}(t)}\mathbf{T}+\mathbf{A}*\mathbf{T},
 \qquad
 \nabla^0\mathbf{g}=\mathbf{A}*\mathbf{g},
 \label{eq:cambio-conexion-continuacion-ricci}
\end{equation}
where, in the second equality, the two summands corresponding to the two
indices of $\mathbf{g}$ are included in the notation $*$. Since $\nabla^0$ is fixed,
\eqref{eq:evolucion-conexion-ricci} reads
\begin{equation}
 \partial_t\mathbf{A}=-\mathbf{g}^{-1}*\nabla^{\mathbf{g}(t)}\operatorname{Ric}.
 \label{eq:evolucion-diferencia-conexiones-ricci}
\end{equation}
The first-order estimate in
Proposition~\ref{prop:estimaciones-derivadas-curvatura-ricci} bounds the
right-hand side and, after integration from $\frac{T_{\displaystyle\max}}{2}$, bounds $\mathbf{A}$.

For higher orders, fix $t_0=T_{\displaystyle\max}/2$ and measure all
infinity norms with $\mathbf{g}(0)$. First specify the conversion of
derivatives. The triangular formula of
Lemma~\ref{lem:local-expression-higher-order}, or its application to the two
connections, takes the form
\begin{equation}
 (\nabla^0)^r\mathbf{T}
 =(\nabla^{\mathbf{g}(t)})^r\mathbf{T}
 +\sum_{j=0}^{r-1}\mathbf{C}_{r,j}*
                   (\nabla^{\mathbf{g}(t)})^j\mathbf{T},
 \qquad r\geq1.
 \label{eq:conversion-orden-r-conexion-fija-ricci}
\end{equation}
Each $\mathbf{C}_{r,j}$ is a finite sum of products of
$\mathbf{A},\nabla^0\mathbf{A},\ldots,
(\nabla^0)^{r-j-1}\mathbf{A}$. For $r=1$ the formula is
\eqref{eq:cambio-conexion-continuacion-ricci}. If it holds for $r$, apply
$\nabla^0$: when the derivative falls on $\mathbf{C}_{r,j}$, it raises
by at most one the order of a derivative of $\mathbf{A}$; when it falls on
$(\nabla^{\mathbf{g}})^j\mathbf{T}$, the first identity in
\eqref{eq:cambio-conexion-continuacion-ricci} gives
$(\nabla^{\mathbf{g}})^{j+1}\mathbf{T}$ and
$\mathbf{A}*(\nabla^{\mathbf{g}})^j\mathbf{T}$. These are precisely the terms
allowed for $r+1$. Contractions are understood as the natural
actions of a difference of connections on every tensor slot;
they introduce no new metric derivatives.

Now suppose that $(\nabla^0)^j\mathbf{A}$ are bounded for
$j\in\{0,\ldots,r-1\}$, with $r\geq1$. Identity
$\nabla^0\mathbf{g}=\mathbf{A}*\mathbf{g}$ shows that derivatives of
$\mathbf{g}$ up to order $r$ are bounded: in the step from order $k$ to
$k+1$, distribute $k$ derivatives between $\mathbf{A}$ and $\mathbf{g}$,
so at most $(\nabla^0)^k\mathbf{A}$ and the already controlled derivatives
of $\mathbf{g}$ occur. For $0\leq k\leq r-1$ all these
factors are available. Differentiating
$\mathbf{g}^{-1}\mathbf{g}=I$ in the same way controls derivatives
of $\mathbf{g}^{-1}$ up to order $r$, also using uniform
equivalence of metrics.

Apply $(\nabla^0)^r$ to
\eqref{eq:evolucion-diferencia-conexiones-ricci}. The Leibniz rule
distributes derivatives between $\mathbf{g}^{-1}$ and
$\nabla^{\mathbf{g}(t)}\operatorname{Ric}$. For the second factor,
\eqref{eq:conversion-orden-r-conexion-fija-ricci} uses only derivatives
of $\mathbf{A}$ up to order $r-1$ and derivatives of Riemann with respect to
$\nabla^{\mathbf{g}(t)}$ up to order $r+1$. By the induction hypothesis
and curvature estimates, all are bounded. Consequently,
\[
 \|\partial_t(\nabla^0)^r\mathbf{A}(t)\|_{L^\infty(M,T^{(1,r+2)}(TM);\mathbf g(0))}
 \leq C_r,\qquad t\in[t_0,T_{\max}).
\]
If $a_r(t):=\|(\nabla^0)^r\mathbf{A}(t)\|_{L^\infty(M,T^{(1,r+2)}(TM);\mathbf g(0))}$,
integration even gives
$a_r(t)\leq a_r(t_0)+C_r(t-t_0)$ and, in particular,
\begin{equation}
 a_r(t)\leq a_r(t_0)+C_r\int_{t_0}^t(1+a_r(\rho))\,d\rho.
 \label{eq:gronwall-conexion-fija-continuacion}
\end{equation}
This proves the inductive step. The case $r=0$ was already obtained by integrating
the evolution of $\mathbf{A}$. All derivatives of
$\mathbf{A}$ with the fixed connection are bounded; the preceding identities then give
those of $\mathbf{g}$, $\mathbf{g}^{-1}$, and $\operatorname{Rm}$.

In particular, conversion of connections gives, for every $r$, a
constant $C_r$ such that
\[
 \|(\nabla^0)^r\partial_t\mathbf{g}(t)\|_{C^0(M,(T^*M)^{\otimes(r+2)};\mathbf{g}(0))}
 =2\|(\nabla^0)^r\operatorname{Ric}_{\mathbf{g}(t)}\|_{C^0(M,(T^*M)^{\otimes(r+2)};\mathbf{g}(0))}
 \leq C_r.
\]
Consequently, if $t_0\leq t_1<t_2<T_{\displaystyle\max}$,
\begin{equation}
 \|\mathbf{g}(t_2)-\mathbf{g}(t_1)\|_{C^r(M,S^2T^*M;\mathbf{g}(0))}
 \leq C_r|t_2-t_1|.
 \label{eq:cauchy-Cr-continuacion-ricci}
\end{equation}
Completeness of each $C^r$ space and compatibility of the limits
produce a smooth metric $\mathbf{g}_*$ such that
\[
 \mathbf{g}(t)\longrightarrow \mathbf{g}_*
 \quad\text{in }C^\infty(M)
 \quad\text{as }t\to T_{\max}^{-}.
\]
Inequality \eqref{eq:equivalencia-uniforme-metricas-ricci} ensures that
$\mathbf{g}_*$ is positive definite.

Theorem~\ref{teo:existencia-unicidad-flujo-ricci}, applied to $\mathbf{g}_*$,
gives a solution on $[T_{\displaystyle\max},T_{\displaystyle\max}+\varepsilon)$ with initial data
$\mathbf{g}_*$. The piecewise definition is continuous, and its first time
derivative agrees on both sides, since it equals $-2\operatorname{Ric}_{\mathbf{g}_*}$.
To check all derivatives, define the differential operators
$\mathcal P_1(\mathbf{g})=-2\operatorname{Ric}_{\mathbf{g}}$ and
\[
 \mathcal P_{k+1}(\mathbf{g})
 :=D\mathcal P_k(\mathbf{g})[-2\operatorname{Ric}_{\mathbf{g}}],
 \qquad k\in\mathbb N.
\]
If $\partial_t^k\mathbf{g}=\mathcal P_k(\mathbf{g})$, the chain
rule gives $\partial_t^{k+1}\mathbf{g}=\mathcal P_{k+1}(\mathbf{g})$.
Thus every time derivative is a smooth expression in the inverse metric
and finitely many spatial derivatives of the metric. Convergence
in all $C^r$ spaces makes these expressions and their spatial
derivatives agree with those of the new solution at $T_{\displaystyle\max}$.
This gives
a smooth extension, contradicting maximality. This proof is
also the geometric application of
Theorem~\ref{teo:criterio-continuacion-parabolico}: the derivative
estimates verify that no coefficient norm of the DeTurck
system diverges while curvature remains bounded. This first proves
that the limit superior is infinite.

If the full limit were not infinite, there would exist $K<\infty$ and
an increasing sequence of times $t_j\to T_{\max}^{-}$ such that
\[
\max_{x\in M}|\operatorname{Rm}_{\mathbf{g}(t_j)}(x)|_{\mathbf{g}(t_j)}\leq K.
\]
Replacing $K$ by $\displaystyle \max\{K,1\}$, we may assume $K>0$.
Apply Proposition~\ref{prop:control-inicial-curvatura-ricci} to the shifted
flow starting at $t_j$. For $j$ large enough that
$T_{\displaystyle\max}-t_j<(2c_mK)^{-1}$, we obtain
\[
\max_{x\in M}|\operatorname{Rm}_{\mathbf{g}(t)}(x)|_{\mathbf{g}(t)}\leq2K,
\qquad t_j\leq t<T_{\max}.
\]
Curvature is also bounded on the compact interval $[0,t_j]$. Thus there is a uniform bound on all of $[0,T_{\displaystyle\max})$, contradicting the
first part of the proof. The conclusion is the stated limit.
\end{proof}

\section{Rescalings and normalized flow}
\label{sec:flujo-ricci-normalizado}

If $c>0$ is constant and $\widehat{\mathbf{g}}=c\mathbf{g}$, the Levi--Civita connections of
$\mathbf{g}$ and $\widehat{\mathbf{g}}$ agree. Consequently,
\begin{equation}
 \operatorname{Rm}^{(1,3)}_{c\mathbf{g}}=\operatorname{Rm}^{(1,3)}_{\mathbf{g}},
\qquad
 \operatorname{Rm}^{(0,4)}_{c\mathbf{g}}=c\,\operatorname{Rm}^{(0,4)}_{\mathbf{g}},
\qquad
 \operatorname{Ric}_{c\mathbf{g}}=\operatorname{Ric}_{\mathbf{g}},
 \qquad
 R_{c\mathbf{g}}=c^{-1}R_{\mathbf{g}},
 \qquad
 d\lambda_{c\mathbf{g}}=c^{\frac{m}{2}}d\lambda_{\mathbf{g}},
\label{eq:reescalamiento-curvatura-ricci}
\end{equation}
as well as
\begin{equation}
 d_{c\mathbf{g}}=c^{1/2}d_{\mathbf{g}},\qquad
 \operatorname{inj}_{c\mathbf{g}}=c^{1/2}\operatorname{inj}_{\mathbf{g}},\qquad
 |\nabla_{c\mathbf{g}}^q\operatorname{Rm}_{c\mathbf{g}}|_{c\mathbf{g}}
 =c^{-1-q/2}|\nabla_{\mathbf{g}}^q\operatorname{Rm}_{\mathbf{g}}|_{\mathbf{g}}.
\label{eq:reescalamiento-derivadas-curvatura-ricci}
\end{equation}
The equality $\operatorname{Ric}_{c\mathbf{g}}=\operatorname{Ric}_{\mathbf{g}}$ is an equality
of covariant tensors of order two; it is important not to replace it by
$c^{-1}\operatorname{Ric}_{\mathbf{g}}$.

Denote by
\begin{equation}
 r(\mathbf{g})=\frac{1}{\operatorname{Vol}(M,\mathbf{g})}
 \int_MR_{\mathbf{g}}\,d\lambda_{\mathbf{g}}
 \label{eq:promedio-curvatura-escalar}
\end{equation}
the average scalar curvature. The volume-normalized Ricci
flow is
\begin{equation}
 \partial_t\mathbf{g}=-2\operatorname{Ric}_{\mathbf{g}}+\frac{2}{m}r(\mathbf{g})\mathbf{g}.
 \label{eq:flujo-ricci-normalizado}
\end{equation}

\begin{proposition}[Volume and curvature of normalized flow]
\label{prop:evolucion-flujo-ricci-normalizado}
Let $M$ be a closed smooth manifold of dimension $m$ and let
$\mathbf{g}(t)$ be a smooth solution of normalized Ricci flow
\eqref{eq:flujo-ricci-normalizado}.
Every solution of \eqref{eq:flujo-ricci-normalizado} satisfies
\begin{align}
 \partial_tg_{ij}
 &=-2\operatorname{Ric}_{ij}+\frac{2}{m}r g_{ij},
 \label{eq:metrica-flujo-normalizado}\\
 \partial_td\lambda_{\mathbf{g}}&=(-R_{\mathbf{g}}+r)d\lambda_{\mathbf{g}},
 \label{eq:medida-flujo-normalizado}\\
 \frac{d}{dt}\operatorname{Vol}(M,\mathbf{g}(t))&=0,
 \label{eq:volumen-flujo-normalizado}\\
 \partial_tR_{\mathbf{g}}
 &=\Delta_{\mathbf{g}} R_{\mathbf{g}}+2|\operatorname{Ric}_{\mathbf{g}}|_{\mathbf{g}}^2-\frac{2}{m}rR_{\mathbf{g}}.
 \label{eq:escalar-flujo-normalizado}
\end{align}
\end{proposition}

\begin{proof}
The first identity is the component expression of
\eqref{eq:flujo-ricci-normalizado}.
The trace of
$\mathbf{h}=-2\operatorname{Ric}_{\mathbf{g}}+\frac{2}{m}r\mathbf{g}$ is
$-2R_{\mathbf{g}}+2r$. By
\eqref{eq:variacion-volumen-metrica} we obtain
\eqref{eq:medida-flujo-normalizado}. Integrating this identity,
\[
 \frac{d}{dt}\operatorname{Vol}(M,\mathbf{g}(t))
 =-\int_MR_{\mathbf{g}}\,d\lambda_{\mathbf{g}}+r\int_Md\lambda_{\mathbf{g}}=0,
\]
where definition \eqref{eq:promedio-curvatura-escalar} was used.
The part $-2\operatorname{Ric}$ of $\mathbf{h}$ gives
$\Delta_{\mathbf{g}} R_{\mathbf{g}}+2|\operatorname{Ric}_{\mathbf{g}}|_{\mathbf{g}}^2$ by
\eqref{eq:evolucion-curvatura-escalar}. For the part
$\frac{2}{m}r\mathbf{g}$, the number $r(t)$ is spatially constant at each time,
so both derivative terms in
\eqref{eq:variacion-curvatura-escalar} vanish, while
\[
 -\left\langle\operatorname{Ric}_{\mathbf{g}},\frac{2}{m}r\mathbf{g}\right\rangle_{\mathbf{g}}
 =-\frac{2}{m}rR_{\mathbf{g}}.
\]
Their sum proves \eqref{eq:escalar-flujo-normalizado}.
\end{proof}

Normalized flow is not a different equation from the local
viewpoint. A spatial rescaling by a positive function of time,
together with a change of the time variable, turns an unnormalized
solution into a normalized one. Indeed, if $\mathbf{g}(t)$ solves
\eqref{eq:flujo-ricci}, set
\[
 c(t):=
 \left(
 \frac{\operatorname{Vol}(M,\mathbf{g}(0))}
 {\operatorname{Vol}(M,\mathbf{g}(t))}
 \right)^{2/m},
 \qquad
 \tau(t):=\int_0^t c(s)\,ds,
 \qquad
 \widehat{\mathbf{g}}(\tau(t)):=c(t)\mathbf{g}(t).
\]
Then $\operatorname{Vol}(\widehat{\mathbf{g}})=\operatorname{Vol}(\mathbf{g}(0))$.
Moreover, viewing $c$ as a function of $\tau$ and $t=t(\tau)$, we have
\[
\frac{dt}{d\tau}=\frac1c,
\qquad
\frac1c\frac{dc}{d\tau}
=\frac{2}{m}r(\widehat{\mathbf{g}}),
\]
since $\partial_t\log\operatorname{Vol}(\mathbf{g})=-r(\mathbf{g})$ and
$r(c\mathbf{g})=c^{-1}r(\mathbf{g})$.
Then, using \eqref{eq:reescalamiento-curvatura-ricci},
\begin{align*}
 \partial_\tau\widehat{\mathbf{g}}
 &=\frac{dc}{d\tau}\mathbf{g}
 +c\frac{dt}{d\tau}\partial_t\mathbf{g}\\
 &=\frac{1}{c}\frac{dc}{d\tau}\widehat{\mathbf{g}}
 -2\operatorname{Ric}_{\mathbf{g}}\\
 &=\frac{2}{m}r(\widehat{\mathbf{g}})\widehat{\mathbf{g}}
 -2\operatorname{Ric}_{\widehat{\mathbf{g}}}.
\end{align*}
This gives precisely \eqref{eq:flujo-ricci-normalizado}.

\section{Einstein metrics, surfaces, and solitons}
\label{sec:ejemplos-flujo-ricci}

The simplest solutions combine rescalings, diffeomorphisms, and
special curvature structures.

\begin{example}[Einstein metrics]
\label{ej:flujo-ricci-metrica-einstein}
Let $M$ be a closed smooth manifold and let $\mathbf{g}_0$ be a Riemannian metric on $M$. Suppose that $\operatorname{Ric}_{\mathbf{g}_0}=\lambda \mathbf{g}_0$, with constant $\lambda$.
On the interval where $1-2\lambda t>0$,
\begin{equation}
 \mathbf{g}(t)=(1-2\lambda t)\mathbf{g}_0
 \label{eq:flujo-ricci-einstein}
\end{equation}
solves Ricci flow. Indeed,
\[
 \partial_t\mathbf{g}=-2\lambda \mathbf{g}_0,
 \qquad
 \operatorname{Ric}_{\mathbf{g}(t)}=\operatorname{Ric}_{\mathbf{g}_0}
 =\lambda \mathbf{g}_0.
\]
If $\lambda>0$, the homothetic solution becomes extinct at
$t=\frac{1}{2\lambda}$; if $\lambda=0$, it is stationary; if $\lambda<0$,
it exists for all $t\geq0$ and expands. The same metric $\mathbf{g}_0$ is a fixed
point of normalized flow, since $r(\mathbf{g}_0)=m\lambda$.
\end{example}

In particular, a metric of constant sectional curvature $K$ is Einstein
with $\lambda=(m-1)K$. Formula
\eqref{eq:flujo-ricci-einstein} describes its evolution completely.

\begin{proposition}[Ricci flow in dimension two]
\label{prop:flujo-ricci-superficies}
On a closed surface,
\begin{align}
 \partial_t\mathbf{g}&=-R\mathbf{g},
 &\partial_tR&=\Delta R+R^2
 \label{eq:flujo-ricci-dimension-dos}
\end{align}
for unnormalized flow. For normalized flow,
\begin{align}
 \partial_t\mathbf{g}&=(r-R)\mathbf{g},
 &\partial_tR&=\Delta R+R(R-r)
 \label{eq:flujo-ricci-normalizado-dimension-dos}
\end{align}
In the latter case,
\begin{equation}
 r=\frac{4\pi\chi(M)}{\operatorname{Vol}(M,\mathbf{g}(0))}
 \label{eq:promedio-escalar-gauss-bonnet}
\end{equation}
is constant in time.
\end{proposition}

\begin{proof}
In dimension two, the algebraic Riemann symmetries imply
\[
 \operatorname{Ric}=K\mathbf{g},
 \qquad R=2K,
\]
and hence $\operatorname{Ric}=\frac{1}{2}R\mathbf{g}$. Substitution into
\eqref{eq:flujo-ricci} and \eqref{eq:evolucion-curvatura-escalar} gives
\[
 \partial_t\mathbf{g}=-R\mathbf{g},
 \qquad
 \partial_tR=\Delta R+2\left|\frac{1}{2}R\mathbf{g}(t)\right|_{\mathbf{g}(t)}^2
 =\Delta R+R^2,
\]
because $|\mathbf{g}(t)|_{\mathbf{g}(t)}^2=2$. The normalized formulas follow in the same way
from \eqref{eq:flujo-ricci-normalizado} and
\eqref{eq:escalar-flujo-normalizado}. If $M$ is orientable,
Theorem~\ref{teo:gauss-bonnet-superficies-apendice} states
\[
 \int_MK\,d\lambda_{\mathbf{g}}=2\pi\chi(M).
\]
If $M$ is nonorientable, apply this equality to the orientable two-sheeted
covering space $\pi\colon\widetilde M\to M$ with the lifted metric. Since
\[
 \int_{\widetilde M}K_{\pi^*\mathbf{g}}\,d\lambda_{\pi^*\mathbf{g}}
 =2\int_MK_{\mathbf{g}}\,d\lambda_{\mathbf{g}},
 \qquad
 \chi(\widetilde M)=2\chi(M),
\]
the same formula holds on $M$. Consequently,
$\displaystyle \int_MR\,d\lambda_{\mathbf{g}}=4\pi\chi(M)$. Volume is constant by
\eqref{eq:volumen-flujo-normalizado}, proving
\eqref{eq:promedio-escalar-gauss-bonnet}.
\end{proof}

\begin{definition}[Ricci soliton]
\label{def:soliton-ricci}
Let $M$ be a smooth manifold with or without boundary. A \emph{Ricci soliton}
on $M$ is a metric $\mathbf{g}_0$ for which there exist a
smooth vector field $\mathbf{X}$ and a constant $\lambda$ such that
\begin{equation}
 \operatorname{Ric}_{\mathbf{g}_0}+\frac{1}{2}\mathcal L_{\mathbf{X}}\mathbf{g}_0=\lambda \mathbf{g}_0.
 \label{eq:ecuacion-soliton-ricci}
\end{equation}
It is called shrinking, steady, or expanding according as
$\lambda>0$, $\lambda=0$, or $\lambda<0$. If $\mathbf{X}=\nabla^{\mathbf{g}_0}f$, the equation
is equivalent to
\begin{equation}
 \operatorname{Ric}_{\mathbf{g}_0}+\nabla^2f=\lambda \mathbf{g}_0
 \label{eq:soliton-gradiente-ricci}
\end{equation}
and it is called a gradient soliton.
\end{definition}

\begin{proposition}[Self-similar solutions]
\label{prop:solucion-autosimilar-soliton-ricci}
Let $M$ be a closed smooth manifold and let
$(\mathbf{g}_0,\mathbf{X},\lambda)$ be a Ricci soliton on $M$. Set
\[
 \sigma(t)=1-2\lambda t
\]
and, while $\sigma(t)>0$, let $\phi_t$ be the family of diffeomorphisms
generated by the time-dependent field
\begin{equation}
 \mathbf{Y}_t=\frac{1}{\sigma(t)}\mathbf{X},
 \qquad \phi_0=\operatorname{id}_M.
 \label{eq:campo-difeomorfismos-soliton}
\end{equation}
Then
\begin{equation}
 \mathbf{g}(t)=\sigma(t)\phi_t^*\mathbf{g}_0
 \label{eq:solucion-autosimilar-soliton}
\end{equation}
is a solution of Ricci flow.
\end{proposition}

\begin{proof}
The pullback differentiation rule and
$\sigma'(t)=-2\lambda$ give
\begin{align*}
 \partial_t\mathbf{g}(t)
 &=-2\lambda\phi_t^*\mathbf{g}_0
 +\sigma(t)\phi_t^*(\mathcal L_{\mathbf{Y}_t}\mathbf{g}_0)\\
 &=\phi_t^*(-2\lambda \mathbf{g}_0+\mathcal L_{\mathbf{X}}\mathbf{g}_0)\\
 &=-2\phi_t^*\operatorname{Ric}_{\mathbf{g}_0},
\end{align*}
where the last equality is
\eqref{eq:ecuacion-soliton-ricci}. By naturality of Ricci and its
invariance under constant rescaling,
\[
 \operatorname{Ric}_{\sigma(t)\phi_t^*\mathbf{g}_0}
 =\operatorname{Ric}_{\phi_t^*\mathbf{g}_0}
 =\phi_t^*\operatorname{Ric}_{\mathbf{g}_0}.
\]
Therefore, $\partial_t\mathbf{g}=-2\operatorname{Ric}_{\mathbf{g}}$.
\end{proof}

Einstein metrics are the solitons with $\mathbf{X}=0$. Solitons with
$\mathbf{X}\neq0$ also incorporate the motion by diffeomorphisms that the geometric
equation does not detect in its principal symbol.
Proposition~\ref{prop:solucion-autosimilar-soliton-ricci} describes the
combination of this motion with rescaling.

\chapter{Compactness, singularities, and Perelman's theory}
\label{cap:compactacion-singularidades-perelman}

When curvature in a Ricci flow grows without bound, the relevant geometry
concentrates in regions whose natural size tends to zero. To observe them,
one chooses a point of large curvature, multiplies the metric by that scale, and
reparametrizes time. If volume does not collapse and curvature derivatives
remain controlled, a subsequence of these rescalings
converges to a complete flow retaining the geometry of the singularity at
unit scale.

This is the first chain of ideas in the chapter: Shi estimates, injectivity
radius, and Hamilton compactness. The second begins with conjugate heat
and Perelman's monotone quantities. The functionals
$\mathcal F$ and $\mathcal W$, the quantity $\mu$, and reduced volume prevent
regions of controlled curvature from losing all their volume. In
dimension three, Hamilton--Ivey pinching turns rescaling
limits into ancient solutions of nonnegative curvature, the
$\kappa$--solutions.

Surgery theory additionally requires a uniform description of regions of
large curvature, the standard solution inserted after cutting a
neck, and an inductive argument preserving pinching, noncollapsing, and
canonical neighborhoods. We first study the analytic part leading
to singularity limits. In the final section, we state with
precise hypotheses the structural surgery and extinction theorems
used to conclude the Poincaré Conjecture; for these results we follow Chapters~12--19 of \cite{morgan2007ricci}.

Perelman introduced these tools in his 2002 and 2003 papers
\cite{Perelman2002,Perelman2003Surgery,Perelman2003Extinction}. The
presentation also draws on the compactness and singularity
estimates in \cite{ChowLuNi2006} and the entropy approach of
\cite{Zhang2010}.

\begin{semblanzaHistorica}{Grisha Perelman and the geometry of singularities}
Perelman's contribution went beyond continuing the flow through
a singularity. He first found quantities distinguishing an
irreversible direction of evolution; he then used noncollapsing and
compactness to prove that magnifying a region of large curvature can
produce only very rigid models. Surgery relies on this rigidity:
one cuts where the geometry is nearly cylindrical and replaces the singular
region by a cap whose behavior can be controlled. In this way,
a singularity ceases to be the point where the argument ends and becomes
an organized part of the evolution.
\end{semblanzaHistorica}

Throughout this chapter, dimension is denoted by $m$.

\section{Pointed smooth convergence}
\label{sec:per-convergencia-apuntada}

A rescaled sequence rarely converges on the entire original manifold.
The appropriate notion fixes a base point and requires smooth convergence only
on increasingly large compact subsets of the limit space. For flows,
spatial identifications must be time-independent; otherwise
the limiting equation would acquire additional Lie derivative terms.

\begin{definition}[Pointed smooth convergence]
\label{def:per-convergencia-suave-apuntada}
A sequence of pointed complete Riemannian manifolds without boundary
$(M_i,\mathbf{g}_i,p_i)$ \emph{converges smoothly in the pointed sense} to
$(M_\infty,\mathbf{g}_\infty,p_\infty)$ if there exist an exhaustion
\[
p_\infty\in U_1\Subset U_2\Subset\cdots\Subset M_\infty,
\qquad
\bigcup_{a=1}^\infty U_a=M_\infty,
\],
a sequence of integers $a_i\to\infty$, and embeddings
$\Phi_i\colon U_{a_i}\to M_i$. For every $a\in\mathbb N$ and every $i$ with
$a\leq a_i$, denote their restrictions by
\[
\Phi_{i,a}\colon U_a\longrightarrow M_i,
\qquad
\Phi_{i,a}(p_\infty)=p_i,
\]
such that
\[
\Phi_{i,a}^*\mathbf{g}_i\longrightarrow \mathbf{g}_\infty
\quad\text{in }C^\infty(U_a).
\]
We also require that the images exhaust the pointed balls of the
approximating spaces: for every $R>0$ there exists $a(R)$ such that, for all sufficiently
large $i$,
\begin{equation}
 \overline{B_{\mathbf{g}_\infty}(p_\infty,R)}\subseteq U_{a(R)},
 \qquad
 B_{\mathbf{g}_i}(p_i,R)\subseteq\Phi_{i,a(R)}(U_{a(R)}).
 \label{eq:per-agotamiento-bolas-convergencia-apuntada}
\end{equation}
The notation $C^\infty(U_a)$ means uniform convergence of every
derivative on each compact subset of $U_a$. Since
$\Phi_{i,a}=\Phi_i|_{U_a}$, the identifications are compatible under
restriction of domains. We work in the connected component of each base
point, and the limit manifold is also taken to be connected and without boundary.
\end{definition}

\begin{definition}[Pointed convergence of flows]
\label{def:per-convergencia-apuntada-flujos}
Let $M_i$ and $M_\infty$ be smooth manifolds without boundary and let
$\mathbf{g}_i(t)$ be complete Ricci flows on $M_i$, defined on a
common interval $I$. We say that
\[
(M_i,\mathbf{g}_i(t),p_i)\longrightarrow
(M_\infty,\mathbf{g}_\infty(t),p_\infty)
\]
smoothly on compact spacetime subsets if the embeddings in the preceding
definition are independent of $t$ and, for every
$J\Subset I$,
\[
\Phi_{i,a}^*\mathbf{g}_i(t)\longrightarrow \mathbf{g}_\infty(t)
\quad\text{in }C^\infty(U_a\times J).
\]
Exhaustion of balls is also required uniformly on compact
time intervals: given $J\Subset I$ and $R>0$, one can choose
$a=a(J,R)$ so that, for sufficiently large $i$ and every $t\in J$,
\begin{equation}
 \overline{B_{\mathbf{g}_\infty(t)}(p_\infty,R)}\subseteq U_a,
 \qquad
 B_{\mathbf{g}_i(t)}(p_i,R)\subseteq\Phi_{i,a}(U_a).
 \label{eq:per-agotamiento-bolas-convergencia-flujos}
\end{equation}
\end{definition}

\begin{lemma}[Passing geometry to the limit]
\label{lem:per-paso-limite-geometria}
Let $M_i$ and $M_\infty$ be smooth manifolds without boundary and suppose that
$(M_i,\mathbf{g}_i(t),p_i)$ converges smoothly in the pointed sense to
$(M_\infty,\mathbf{g}_\infty(t),p_\infty)$. Then, on every identified compact
subset, the Levi--Civita connections, the
$\operatorname{Rm}$ and $\operatorname{Ric}$ tensors, scalar curvatures, and
all their covariant derivatives converge. Convergence also commutes with
tensor products, contractions, and traces. If the $\mathbf{g}_i(t)$
are Ricci flows and convergence is
in spacetime, the limit satisfies
$\partial_t\mathbf{g}_\infty=-2\operatorname{Ric}_{\mathbf{g}_\infty}$.
\end{lemma}

\begin{proof}
In a fixed chart, Christoffel symbols are smooth rational
functions of $g_{ij}$, $g^{ij}$, and $\partial g_{ij}$. Curvature and its
derivatives are obtained by iterating derivatives and products of these coefficients.
$C^\infty$ convergence allows passage to the limit in all these
expressions. Tensor operations are continuous in each
$C^k$ space on compact subsets. Finally, $C^\infty$ convergence
in the time variable allows passage to the limit in the flow
equation.
\end{proof}

The injectivity radius prevents a curvature bound from concealing collapse.

\begin{lemma}[Volume, curvature, and injectivity radius]
\label{lem:per-volumen-curvatura-inyectividad}
For every $m\in\mathbb N$ and $v>0$ there exists $c=c(m,v)>0$ with the following
property. If $(N,\mathbf{h})$ is a complete Riemannian manifold of dimension $m$
without boundary,
\[
 |\operatorname{Rm}_{\mathbf{h}}|_{\mathbf{h}}\leq r^{-2}
 \quad\hbox{on }B_{\mathbf{h}}(x,r),
 \qquad
 \operatorname{Vol}_{\mathbf{h}}B_{\mathbf{h}}(x,r)\geq vr^m,
\]
then $\operatorname{inj}_{\mathbf{h}}(x)\geq cr$.
\end{lemma}

The preceding estimate is local: it involves only the geometry of the ball
$B_{\mathbf{h}}(x,r)$ and completeness of the metric. Its proof, based
on the short-loop lemma and comparison of the exponential
Jacobian, can be found in Chapter~1, in the section
\emph{Local volume and the injectivity radius}, of
\cite{morgan2007ricci}.

\begin{lemma}[Local propagation of injectivity radius]
\label{lem:per-propagacion-local-inyectividad}
Given $m\in\mathbb N$, $R>0$, $\iota_0>0$, and $0\leq K<\infty$, there exists
$\iota=\iota(m,R,\iota_0,K)>0$ with the following property. If
$(N,\mathbf{h},p)$ is a complete Riemannian manifold of dimension $m$ without boundary,
\[
 \operatorname{inj}_{\mathbf{h}}(p)\geq\iota_0,
 \qquad
 |\operatorname{Rm}_{\mathbf{h}}|_{\mathbf{h}}\leq K
 \quad\text{on }B_{\mathbf{h}}(p,2R+2),
\]
then
\[
 \operatorname{inj}_{\mathbf{h}}(x)\geq\iota
 \qquad\text{for every }x\in B_{\mathbf{h}}(p,R).
\]
\end{lemma}

\begin{proof}
After increasing $K$ by a constant depending only on dimension,
we may assume that
\[
 \operatorname{Sec}_{\mathbf{h}}\leq K,
 \qquad
 \operatorname{Ric}_{\mathbf{h}}\geq-(m-1)K\mathbf{h}
\]
on the indicated ball. Choose
\[
 0<r_0<\min\left\{\frac{\iota_0}{2},1,
                     \frac{\pi}{4\sqrt{K+1}}\right\}.
\]
The exponential map at $p$ is a diffeomorphism on the tangent ball of
radius $2r_0$. The local Jacobi-field argument used in the proof
of Theorem~\ref{teo: equivalencias de geometria acotada}, applied on this
ball, gives a constant $c_0=c_0(m,K,r_0)>0$ such that
\[
 \operatorname{Vol}_{\mathbf{h}}B_{\mathbf{h}}(p,r_0)\geq c_0r_0^m.
\]

Let $x\in B_{\mathbf{h}}(p,R)$. Since
$B_{\mathbf{h}}(p,r_0)\subseteq B_{\mathbf{h}}(x,R+r_0)$, we have
\[
 \operatorname{Vol}_{\mathbf{h}}B_{\mathbf{h}}(x,R+r_0)
 \geq c_0r_0^m.
\]
Denote by $V_{-K}(s)$ the volume of a ball of radius $s$ in the simply
connected space of constant sectional curvature $-K$. Relative
Bishop--Gromov comparison, valid because
$\operatorname{Ric}_{\mathbf{h}}\geq-(m-1)K\mathbf{h}$, gives
\[
 \frac{\operatorname{Vol}_{\mathbf{h}}B_{\mathbf{h}}(x,\rho)}
      {V_{-K}(\rho)}
 \geq
 \frac{\operatorname{Vol}_{\mathbf{h}}B_{\mathbf{h}}(x,R+r_0)}
      {V_{-K}(R+r_0)}
\]
for $0<\rho\leq R+r_0$. For this form of comparison, see
\cite[Chapter~1]{morgan2007ricci}. Take
\[
 \rho:=\min\{1,(K+1)^{-1/2},R+r_0\}.
\]
The preceding estimate implies
\[
 \operatorname{Vol}_{\mathbf{h}}B_{\mathbf{h}}(x,\rho)
 \geq v\rho^m
\]
with $v=v(m,R,\iota_0,K)>0$. Moreover,
$B_{\mathbf{h}}(x,\rho)\subseteq B_{\mathbf{h}}(p,R+1)$ and
$|\operatorname{Rm}_{\mathbf{h}}|\leq K\leq\rho^{-2}$ on this ball.
Lemma~\ref{lem:per-volumen-curvatura-inyectividad} gives
\[
 \operatorname{inj}_{\mathbf{h}}(x)\geq c(m,v)\rho.
\]
The constant on the right is independent of $x$ and the manifold.
\end{proof}

The version of Shi's estimates proved in the preceding chapter
was formulated for closed manifolds. For compactness we will need the
complete noncompact version.

\begin{theorem}[Shi estimates for complete flows]
\label{teo:per-shi-completo}
Let $(M,\mathbf{g}(t))$ (of dimension $m$), $t\in[s_0,s_1]$, be a complete Ricci flow and
suppose that
\[
 \sup_{t\in[s_0,s_1]}
 \|\operatorname{Rm}_{\mathbf{g}(t)}\|_{L^\infty(M,\mathbf{g}(t))}
 \leq K<\infty.
\]
For every $q\in\mathbb N$ and every $\delta\in(0,s_1-s_0]$ there exists a
constant $C_q=C_q(m,q,K,s_1-s_0,\delta)>0$ such that
\[
 \sup_{t\in[s_0+\delta,s_1]}
 \|\nabla_{\mathbf{g}(t)}^q\operatorname{Rm}_{\mathbf{g}(t)}\|_{L^\infty(M,\mathbf{g}(t))}
 \leq C_q.
\]
In particular, if $J\Subset I$ and there exists $\displaystyle s<\min J$ such that
$\displaystyle [s,\displaystyle\max J]\subseteq I$, a uniform curvature bound on this
interval gives uniform bounds for all its covariant derivatives
on $M\times J$.
\end{theorem}

We will use this theorem in its complete version, proved through
local estimates in \cite[Chapter~3, Section~7]{morgan2007ricci}.
The Bernstein quantities are those of
Proposition~\ref{prop:estimaciones-derivadas-curvatura-ricci};
localization allows the maximum principle to be applied on compact domains
before their radius is sent to infinity. The time margin $\delta>0$
is essential: the constant may diverge as $\delta\to0^{+}$.
In particular, a curvature bound with no initial bounds for its derivatives
only gives smooth convergence at times later than the initial time.

We separate from the parabolic argument the compactness theorem for a
single metric.

\begin{theorem}[Cheeger--Gromov pointed smooth compactness]
\label{teo:per-compactacion-estatica-cheeger-gromov}
Let $(M_i^m,\mathbf{h}_i,p_i)$ be a sequence of complete Riemannian
manifolds without boundary. Suppose that there exists $\iota_0>0$ such that
\[
 \operatorname{inj}_{\mathbf{h}_i}(p_i)\geq\iota_0
 \qquad\text{for every }i,
\]
and that, for every $R>0$ and $q\in\mathbb N_0$, there exists
$C_{q,R}<\infty$ such that
\[
 \sup_{i\in\mathbb N}
 \|\nabla_{\mathbf{h}_i}^q\operatorname{Rm}_{\mathbf{h}_i}\|_
 {L^\infty(B_{\mathbf{h}_i}(p_i,R),\mathbf{h}_i)}\leq C_{q,R}.
\]
Then there exist a complete Riemannian manifold without boundary
$(M_\infty,\mathbf{h}_\infty,p_\infty)$ and a subsequence converging
smoothly to it in the pointed sense. The convergence embeddings can
be chosen compatibly on an exhaustion
\[
 p_\infty\in U_1\Subset U_2\Subset\cdots\Subset M_\infty
\]
satisfying the ball-exhaustion property
\eqref{eq:per-agotamiento-bolas-convergencia-apuntada}.
\end{theorem}

\begin{proof}
Replace each $M_i$ by the connected component of $p_i$. Fix
$R>0$. The hypothesis with radius $2R+4$ and
Lemma~\ref{lem:per-propagacion-local-inyectividad} give
$\iota_R>0$ such that
\[
 \operatorname{inj}_{\mathbf{h}_i}(x)\geq\iota_R
 \qquad
 \bigl(x\in B_{\mathbf{h}_i}(p_i,R+1)\bigr)
\]
for every $i$. Bounds for curvature derivatives on the ball of
radius $R+2$, together with the local part of the proof of
Theorem~\ref{teo: equivalencias de geometria acotada}, allow a
radius $\rho_R>0$ to be chosen independently of $i$ and the center
$x\in B_{\mathbf{h}_i}(p_i,R+1)$, such that normal coordinates of radius
$4\rho_R$ are defined and the coefficients of $\mathbf{h}_i$,
its inverse, and all their jets are uniformly bounded in them. The same
bounds hold for transitions between two such charts.

Choose in $B_{\mathbf{h}_i}(p_i,R)$ a maximal net of points separated
by $\rho_R$. The balls of radius $\rho_R$ centered at net points cover the
pointed ball, while the concentric balls of radius
$\rho_R/2$ are disjoint. The Bishop--Gromov argument used in the
preceding lemma gives a uniform lower bound for the volume of
each small ball. The Ricci lower bound gives,
by the same comparison, a uniform upper bound for the volume of
$B_{\mathbf{h}_i}(p_i,R+1)$. Consequently, the number of net
points is bounded by a constant $N_R$ independent of $i$. Thus,
for every $R$, we obtain a finite atlas with uniformly bounded
number of charts, common radius, and uniform bounds for all jets of the metrics
and transitions.

To turn these uniform atlases into a limit manifold, apply
the geometric convergence theorem of
\cite[Theorem~5.6]{morgan2007ricci}. In the complete case under
consideration, its hypotheses are: a nondegenerate base point, upper
volume bounds on fixed-radius balls, bounds for all curvature
derivatives on these balls, and lower volume bounds at small
radii. The injectivity bound verifies the first; Bishop--Gromov
comparison and the preceding paragraph's construction verify the second
and fourth; the third is our hypothesis. Completeness also ensures
that closed balls in every approximating manifold are compact.

The theorem gives a smooth Hausdorff manifold $M_\infty$, a
metric $\mathbf{h}_\infty$, relatively compact open subsets
$U_1\Subset U_2\Subset\cdots$, and, after choosing a subsequence,
embeddings $\Phi_i\colon U_{a_i}\to M_i$, with $a_i\to\infty$.
For every $a\in\mathbb N$, the restrictions
$\Phi_{i,a}=\Phi_i|_{U_a}$ satisfy
\[
 \Phi_{i,a}^*\mathbf{h}_i\longrightarrow\mathbf{h}_\infty
 \quad\text{in }C^\infty(U_a).
\]
The topological part of that theorem constructs the transition domains
and embeddings simultaneously; it therefore also preserves
exhaustion of balls in the approximating manifolds. This is the complete form of
Cheeger--Gromov compactness stated in
\cite[Theorem~5.9]{morgan2007ricci}.

In the cited construction, the atlases of radius $R+1$ cover the ball
$B_{\mathbf{h}_i}(p_i,R+1)$ before the charts are restricted. Reindexing
the exhaustion gives, for every $R\in\mathbb N$ and all sufficiently large $i$,
\begin{equation}
 \overline{B_{\mathbf{h}_\infty}(p_\infty,R)}\subseteq U_{R+1},
 \qquad
 B_{\mathbf{h}_i}(p_i,R)\subseteq\Phi_{i,R+1}(U_{R+1}).
 \label{eq:per-agotamiento-bolas-compactacion-estatica}
\end{equation}
The second inclusion follows by retaining all net charts
meeting the ball of radius $R$; the unit margin allows their domains to be
restricted without leaving points of that ball outside the image. This verifies
\eqref{eq:per-agotamiento-bolas-convergencia-apuntada}. The first inclusion
also shows that every closed bounded ball in
$(M_\infty,\mathbf{h}_\infty)$ is compact. The Hopf--Rinow theorem implies
that $\mathbf{h}_\infty$ is complete.
\end{proof}

The next estimate compares different time slices of a
flow.

\begin{lemma}[Metric and distance distortion]
\label{lem:per-distorsion-metrica-distancias}
Let $(M,\mathbf{g}(t))$, $t\in J$, be a Ricci flow and suppose that
\[
 |\operatorname{Ric}_{\mathbf{g}(t)}|_{\mathbf{g}(t)}\leq K
 \qquad\text{on }M\times J.
\]
For $s,t\in J$ we have
\begin{equation}
 e^{-2K|t-s|}\mathbf{g}(s)
 \leq\mathbf{g}(t)
 \leq e^{2K|t-s|}\mathbf{g}(s),
 \label{eq:per-distorsion-metricas}
\end{equation}
and, if the distances are finite,
\begin{equation}
 e^{-K|t-s|}d_{\mathbf{g}(s)}(x,y)
 \leq d_{\mathbf{g}(t)}(x,y)
 \leq e^{K|t-s|}d_{\mathbf{g}(s)}(x,y).
 \label{eq:per-distorsion-distancias}
\end{equation}
In particular,
\begin{equation}
 B_{\mathbf{g}(s)}\bigl(x,e^{-K|t-s|}r\bigr)
 \subseteq B_{\mathbf{g}(t)}(x,r)
 \subseteq B_{\mathbf{g}(s)}\bigl(x,e^{K|t-s|}r\bigr).
 \label{eq:per-inclusiones-bolas-tiempo}
\end{equation}
If one metric $\mathbf{g}(s)$ is complete, every metric in the
family is complete.
\end{lemma}

\begin{proof}
Fix $v\in T_xM\setminus\{0\}$. Since
$\partial_t\mathbf{g}=-2\operatorname{Ric}_{\mathbf{g}(t)}$,
\[
 \left|\frac d{dt}\log\mathbf{g}(t)(v,v)\right|
 =2\frac{|\operatorname{Ric}_{\mathbf{g}(t)}(v,v)|}
          {\mathbf{g}(t)(v,v)}
 \leq2K.
\]
Integration between $s$ and $t$ gives
\[
 -2K|t-s|
 \leq\log\frac{\mathbf g(t)_x(v,v)}{\mathbf g(s)_x(v,v)}
 \leq2K|t-s|.
\]
The exponential function is increasing. Applying it and multiplying
by $\mathbf g(s)_x(v,v)>0$ yields
\[
 e^{-2K|t-s|}\mathbf g(s)_x(v,v)
 \leq\mathbf g(t)_x(v,v)
 \leq e^{2K|t-s|}\mathbf g(s)_x(v,v).
\]
The inequality also holds for $v=0$. Since the
constants are independent of $x$ and $v$, this is precisely
\eqref{eq:per-distorsion-metricas} with the interpretation of
Definition~\ref{def:orden-formas-metricas}.

Taking square roots gives
\[
 e^{-K|t-s|}|v|_{\mathbf g(s)}
 \leq|v|_{\mathbf g(t)}
 \leq e^{K|t-s|}|v|_{\mathbf g(s)}.
\]
For a piecewise smooth curve $\gamma\colon[a,b]\to M$, apply
this bound to $v=\dot\gamma(r)$ on each smooth segment and integrate
with respect to $r$. The identity
$L_{\mathbf g(t)}(\gamma)=\int_a^b|\dot\gamma(r)|_{\mathbf g(t)}\,dr$
gives
\[
 e^{-K|t-s|}L_{\mathbf{g}(s)}(\gamma)
 \leq L_{\mathbf{g}(t)}(\gamma)
 \leq e^{K|t-s|}L_{\mathbf{g}(s)}(\gamma).
\]
Taking the infimum over curves joining $x$ to $y$ yields
\eqref{eq:per-distorsion-distancias}; the ball inclusions are an
immediate consequence.

Suppose that $\mathbf{g}(s)$ is complete. A Cauchy sequence for
$d_{\mathbf{g}(t)}$ is Cauchy for $d_{\mathbf{g}(s)}$ by
\eqref{eq:per-distorsion-distancias} and therefore converges with respect to
$d_{\mathbf{g}(s)}$. The two distances induce the same topology, and the
second inequality in \eqref{eq:per-distorsion-distancias} gives
convergence with respect to $d_{\mathbf{g}(t)}$. Thus $\mathbf{g}(t)$ is
complete.
\end{proof}

We also need to turn Shi's covariant bounds into coefficient bounds
in time-independent coordinates.

\begin{lemma}[Control with respect to a reference connection]
\label{lem:per-control-conexion-referencia}
Let $J\subseteq\mathbb R$ be a compact interval containing $0$ and let
$\mathbf{g}_i(t)$ be Ricci flows defined on open subsets $V_i$, with
$\mathbf{h}_i:=\mathbf{g}_i(0)$. Suppose that the metrics $\mathbf{h}_i$
have uniform bounds for all jets of their coefficients and inverses in
fixed charts over $V_i$. Also suppose that, for every
$q\in\mathbb N_0$, there exists $C_q$ such that
\[
 \sup_{i\in\mathbb N}\ \sup_{t\in J}
 \|\nabla_{\mathbf{g}_i(t)}^q\operatorname{Rm}_{\mathbf{g}_i(t)}\|_
 {L^\infty(V_i,\mathbf{g}_i(t))}\leq C_q.
\]
Then, in each coordinate subdomain whose closure lies within the
fixed charts, all coefficients of $\mathbf{g}_i(t)$ and
$\mathbf{g}_i(t)^{-1}$ and all their spatial and time derivatives
are uniformly bounded on $J$.
\end{lemma}

\begin{proof}
The order-zero curvature bound gives
$|\operatorname{Ric}_{\mathbf g_i(t)}|_{\mathbf g_i(t)}
\leq K:=\sqrt m\,C_0$ for all indices and times.
Since $0\in J$, Lemma~\ref{lem:per-distorsion-metrica-distancias}
implies, with $Q=e^{2K\operatorname{diam}J}$,
\[
 Q^{-1}(\mathbf h_i)_x(v,v)
 \leq\mathbf g_i(t)_x(v,v)
 \leq Q(\mathbf h_i)_x(v,v)
 \qquad(i\in\mathbb N,\ t\in J,\ x\in V_i,\ v\in T_xV_i).
\]
By Proposition~\ref{prop:comparacion-metricas-duales-tensores-volumen},
norms of the same tensor of fixed type are comparable with
constants independent of $i,x,t$. This allows the
metric used in a tensor bound to be changed; to change
covariant derivatives, we study the difference of connections.
Let $\nabla^i(t)$ and $\nabla^{i,0}$ be their Levi--Civita connections and
set
\[
 \mathbf{A}_i(t):=\nabla^i(t)-\nabla^{i,0}.
\]
Then $\mathbf{A}_i(0)=0$, and identity
\eqref{eq:evolucion-conexion-ricci} takes the form
\begin{equation}
 \partial_t\mathbf{A}_i
 =-\mathbf{g}_i(t)^{-1}*\nabla^i(t)\operatorname{Ric}_{\mathbf{g}_i(t)}.
 \label{eq:per-evolucion-diferencia-conexiones}
\end{equation}
Moreover,
\begin{equation}
 \nabla^{i,0}\mathbf{T}
 =\nabla^i(t)\mathbf{T}+\mathbf{A}_i*\mathbf{T},
 \qquad
 \nabla^{i,0}\mathbf{g}_i(t)=\mathbf{A}_i*\mathbf{g}_i(t).
 \label{eq:per-cambio-conexion-referencia}
\end{equation}

Write $D_i=\nabla^{i,0}$ and $\nabla_i=\nabla^i(t)$.
The difference of connections acts on every tensor slot as in
Chapter~\ref{cap:sobolev-haces}. Apply the formula already proved in
\eqref{eq:conversion-orden-r-conexion-fija-ricci}: for every integer
$r\geq1$,
\begin{equation}
 D_i^r\mathbf{T}
 =\nabla_i^r\mathbf{T}
  +\sum_{j=0}^{r-1}\mathbf{C}_{i,r,j}*\nabla_i^j\mathbf{T}.
 \label{eq:per-conversion-derivadas-referencia}
\end{equation}
Each coefficient $\mathbf{C}_{i,r,j}$ is a finite sum of products of
$D_i^b\mathbf{A}_i$ with $b\in\{0,\ldots,r-j-1\}$; the number of
terms depends only on $r$, dimension, and the type of $\mathbf{T}$.
The Leibniz rule of Proposition~\ref{prop: leibniz operador *}
and the contraction estimates in the chapter on bundles bound
these products by the norms of their factors.

The induction starts at $r=0$. The evolution of $\mathbf{A}_i$ and the bound
for $\nabla_i\operatorname{Ric}$ give
$\|\partial_t\mathbf{A}_i(t)\|_{L^\infty(V_i,T^{(1,2)}(TV_i);\mathbf h_i)}\leq C$.
Integrating from $0$ and using $\mathbf{A}_i(0)=0$ yields
$\|\mathbf{A}_i(t)\|_{L^\infty(V_i,T^{(1,2)}(TV_i);\mathbf h_i)}\leq C|J|$.

Now let $r\geq1$ and suppose that the derivatives $D_i^b\mathbf{A}_i$ are bounded uniformly in $i$ and
$t\in J$ for
$b\in\{0,\ldots,r-1\}$. Define
$\mathbf{F}_i=\partial_t\mathbf{A}_i$.
Equation \eqref{eq:per-evolucion-diferencia-conexiones} expresses
$\mathbf{F}_i$ as contractions of $\mathbf{g}_i^{-1}$ and
$\nabla_i\operatorname{Ric}$. Since $\nabla_i\mathbf{g}_i=0$,
for each $j\in\{0,\ldots,r\}$, $\nabla_i^j\mathbf{F}_i$ contains only
Ricci derivatives of order $j+1$ contracted with the inverse
metric. All are bounded by the curvature hypothesis.
Applying \eqref{eq:per-conversion-derivadas-referencia} to
$\mathbf{T}=\mathbf{F}_i$, the additional coefficients contain only
derivatives of $\mathbf{A}_i$ of order at most $r-1$. By the induction
hypothesis,
\[
 \|\partial_tD_i^r\mathbf{A}_i(t)\|_{L^\infty(V_i,T^{(1,r+2)}(TV_i);\mathbf h_i)}
 =\|D_i^r\mathbf{F}_i(t)\|_{L^\infty(V_i,T^{(1,r+2)}(TV_i);\mathbf h_i)}\leq C_r,
 \qquad i\in\mathbb N,\quad t\in J.
\]
Since $D_i^r\mathbf{A}_i(0)=0$, integration toward either
endpoint of $J$ gives
$\|D_i^r\mathbf{A}_i(t)\|_{L^\infty(V_i,T^{(1,r+2)}(TV_i);\mathbf h_i)}\leq C_r|J|$.
This proves the step from $r-1$ to $r$ without introducing a derivative of
$\mathbf{A}_i$ of higher order than the one controlled.

Now apply \eqref{eq:per-conversion-derivadas-referencia} to
$\mathbf{T}=\mathbf{g}_i$ and $\mathbf{T}=\mathbf{g}_i^{-1}$.
All their positive-order derivatives with respect to $\nabla_i$ vanish.
Thus $D_i^r\mathbf{g}_i$ and $D_i^r\mathbf{g}_i^{-1}$ are expressed
in terms of $\mathbf{C}_{i,r,0}$ and already bounded order-zero tensors.
The local expression for higher covariant derivatives in
Lemma~\ref{lem:local-expression-higher-order}, applied to the fixed
connection $D_i$, converts these bounds into coordinate derivative bounds.
Its coefficients are uniform because the jets of $\mathbf{h}_i$ and
$\mathbf{h}_i^{-1}$ are uniform by hypothesis.

For time derivatives, use the coordinate expression for Ricci:
it is a finite sum of products of coefficients of $\mathbf{g}_i^{-1}$,
first derivatives of $\mathbf{g}_i$, and second derivatives of
$\mathbf{g}_i$. Suppose that all spatial derivatives of
$\partial_t^b\mathbf{g}_i$ and $\partial_t^b\mathbf{g}_i^{-1}$ are controlled for
$b\in\{0,\ldots,a\}$. Applying $\partial_t^a$ to
$\partial_t\mathbf{g}_i=-2\operatorname{Ric}_{\mathbf{g}_i}$, each
product of $d$ factors is differentiated by
\[
 \partial_t^a(F_1\cdots F_d)
 =\sum_{\substack{(b_1,\ldots,b_d)\in\mathbb N_0^d\\
                   b_1+\cdots+b_d=a}}
   \frac{a!}{b_1!\cdots b_d!}
   (\partial_t^{b_1}F_1)\cdots(\partial_t^{b_d}F_d).
\]
All time indices on the right-hand side are at most $a$.
The induction hypothesis therefore bounds all spatial derivatives of
$\partial_t^{a+1}\mathbf{g}_i$. Differentiating the matrix identity
$\mathbf{g}_i^{-1}\mathbf{g}_i=I$ then gives
\[
 \partial_t^{a+1}\mathbf{g}_i^{-1}
 =-\sum_{b=0}^{a}\binom{a+1}{b}
   (\partial_t^b\mathbf{g}_i^{-1})
   (\partial_t^{a+1-b}\mathbf{g}_i)\mathbf{g}_i^{-1}.
\]
All factors here are also controlled, even after spatial
derivatives are applied. The initial case $a=0$ is the spatial control
already proved. Bounds for all mixed orders follow.
\end{proof}

\begin{theorem}[Hamilton compactness]
\label{teo:per-compactacion-hamilton}
Let $M_i$ be smooth manifolds without boundary of dimension $m$ and let
$(M_i,\mathbf{g}_i(t),p_i)$ be complete Ricci flows defined on a
common interval $I$, open on the left, with $0\in I$ and
$\displaystyle \inf I<0$. The interval $I$ may contain its right endpoint, in particular
$I=(a,0]$ with $-\infty\leq a<0$. In this statement, $J\Subset I$
denotes a compact interval contained in $I$; it may include the right
endpoint if this belongs to $I$. Suppose that
\[
 \operatorname{inj}_{\mathbf{g}_i(0)}(p_i)\geq\iota_0>0.
\]
Also suppose that, for every $J\Subset I$, there exists $C(J)<\infty$ such that
\[
 \sup_{i\in\mathbb N}\ \sup_{t\in J}
 \|\operatorname{Rm}_{\mathbf{g}_i(t)}\|_{L^\infty(M_i,\mathbf{g}_i(t))}
 \leq C(J).
\]
Then a subsequence converges smoothly, in the pointed sense and on
compact spacetime subsets, to a complete Ricci flow
\[
 (M_\infty,\mathbf{g}_\infty(t),p_\infty),
 \qquad t\in I,
\]
where $M_\infty$ is a smooth manifold without boundary and
$\mathbf{g}_\infty(t)$ is complete for every $t$.
\end{theorem}

\begin{proof}
Set $\mathbf{h}_i:=\mathbf{g}_i(0)$. First extract a limit of
the time slices at $t=0$. Fix $q\in\mathbb N_0$ and choose
$s<0$ with $[s,0]\subseteq I$. The hypothesis gives a uniform curvature
bound on $M_i\times[s,0]$. If $q\geq1$, apply
Theorem~\ref{teo:per-shi-completo} to the shifted flow starting at
$s$ and evaluate at $t=0$; for $q=0$, use the hypothesis directly.
We obtain a constant $C_q$ independent of $i$ such that
\[
 \|\nabla_{\mathbf{h}_i}^q\operatorname{Rm}_{\mathbf{h}_i}\|_
 {L^\infty(M_i,\mathbf{h}_i)}\leq C_q,
 \qquad i\in\mathbb N.
\]
Together with the injectivity radius bound at $p_i$, these estimates
verify the hypotheses of
Theorem~\ref{teo:per-compactacion-estatica-cheeger-gromov}. After
passing to a subsequence, there exist a complete Riemannian manifold
$(M_\infty,\mathbf{h}_\infty,p_\infty)$, an exhaustion
\[
 p_\infty\in U_1\Subset U_2\Subset\cdots\Subset M_\infty
\],
and compatible embeddings
\[
 \Phi_{i,a}\colon U_a\longrightarrow M_i,
 \qquad \Phi_{i,a}(p_\infty)=p_i,
\]
such that
\begin{equation}
 \Phi_{i,a}^*\mathbf{h}_i\longrightarrow\mathbf{h}_\infty
 \quad\text{in }C^\infty(U_a).
 \label{eq:per-convergencia-rebanada-cero-hamilton}
\end{equation}

Now fix $a\in\mathbb N$, a compact set $K\Subset U_a$, and a
compact interval $J\Subset I$. Let $\widetilde J\Subset I$ be a compact
interval containing $J\cup\{0\}$. Choose $\displaystyle s<\min\widetilde J$ so
that $\displaystyle [s,\displaystyle\max\widetilde J]\subseteq I$. The curvature bound on this
interval and Theorem~\ref{teo:per-shi-completo} give, for every
$q\in\mathbb N_0$, constants $C_{q,\widetilde J}$ such that
\begin{equation}
 \sup_{i\in\mathbb N}\ \sup_{t\in\widetilde J}
 \|\nabla_{\mathbf{g}_i(t)}^q\operatorname{Rm}_{\mathbf{g}_i(t)}\|_
 {L^\infty(M_i,\mathbf{g}_i(t))}\leq C_{q,\widetilde J}.
 \label{eq:per-cotas-shi-hamilton}
\end{equation}
Consider on $U_a$ the metrics
\[
 \widehat{\mathbf{g}}_{i,a}(t):=\Phi_{i,a}^*\mathbf{g}_i(t),
 \qquad
 \widehat{\mathbf{h}}_{i,a}:=\widehat{\mathbf{g}}_{i,a}(0).
\]
By \eqref{eq:per-convergencia-rebanada-cero-hamilton}, in any
finite atlas covering $K$, the coefficients of
$\widehat{\mathbf{h}}_{i,a}$ and their inverses have uniform bounds of
all orders. The estimates
\eqref{eq:per-cotas-shi-hamilton} are preserved under pullback.
Lemma~\ref{lem:per-control-conexion-referencia}, applied on
$\widetilde J$, then shows that on $K\times J$ all coefficients of
$\widehat{\mathbf{g}}_{i,a}(t)$ and all their spatial and time derivatives
are uniformly bounded. Metric distortion also ensures
a uniform lower bound for these positive definite tensors.

Exhaust each $U_a$ by compact sets and $I$ by compact
intervals. Arzelà--Ascoli, applied successively to all
coefficient derivatives in finitely many charts, and a
diagonal argument give a subsequence for which
\begin{equation}
 \Phi_{i,a}^*\mathbf{g}_i(t)
 \longrightarrow\mathbf{g}_\infty(t)
 \quad\text{in }C^\infty(K\times J)
 \label{eq:per-convergencia-flujos-hamilton}
\end{equation}
for every $a$, every $K\Subset U_a$, and every $J\Subset I$. Compatibility
of the embeddings $\Phi_{i,a}$ makes the limits obtained on $U_a$ and
$U_{a+1}$ agree on $U_a$; they therefore define a global family of
smooth metrics $\mathbf{g}_\infty(t)$ on $M_\infty$. At $t=0$, identity
\eqref{eq:per-convergencia-rebanada-cero-hamilton} gives
$\mathbf{g}_\infty(0)=\mathbf{h}_\infty$.
Lemma~\ref{lem:per-paso-limite-geometria}, applied to
\eqref{eq:per-convergencia-flujos-hamilton}, shows that
\[
 \partial_t\mathbf{g}_\infty(t)
 =-2\operatorname{Ric}_{\mathbf{g}_\infty(t)}.
\]
Let us also verify uniform exhaustion of balls. Fix
$J\Subset I$ and $R>0$. Let $\widehat J$ be the compact interval joining
$J$ to $0$. The curvature bound on $M_i\times\widehat J$ bounds
$|\operatorname{Ric}|$ by a constant $K_J$ independent of $i$.
By
\eqref{eq:per-distorsion-distancias}, for every $t\in J$,
\[
 B_{\mathbf{g}_i(t)}(p_i,R)
 \subseteq B_{\mathbf{h}_i}(p_i,e^{K_J|t|}R).
\]
The exhaustion property of the static limit allows a common
$U_a$ to be chosen whose image contains the balls on the right for all
$t\in J$. Passing to the limit, the analogous inequality between
$\mathbf{g}_\infty(t)$ and $\mathbf{h}_\infty$ shows that the same $U_a$
contains $\overline{B_{\mathbf{g}_\infty(t)}(p_\infty,R)}$. Thus
\eqref{eq:per-agotamiento-bolas-convergencia-flujos} holds, establishing pointed smooth
convergence on compact spacetime subsets.

It remains to prove completeness of each slice. Fix $t\in I$ and choose
a compact interval $J\Subset I$ containing $0$ and $t$. The curvature
bound gives $K_J<\infty$ such that
\[
 |\operatorname{Ric}_{\mathbf{g}_i(\rho)}|_{\mathbf{g}_i(\rho)}
 \leq K_J
 \qquad ((x,\rho)\in M_i\times J)
\]
for every $i$. By
Lemma~\ref{lem:per-distorsion-metrica-distancias},
\[
 e^{-2K_J|t|}\mathbf{h}_i
 \leq\mathbf{g}_i(t)
 \leq e^{2K_J|t|}\mathbf{h}_i.
\]
Pull back by $\Phi_{i,a}$ and pass to the limit. Since the
$U_a$ cover $M_\infty$, we obtain
\[
 e^{-2K_J|t|}\mathbf{h}_\infty
 \leq\mathbf{g}_\infty(t)
 \leq e^{2K_J|t|}\mathbf{h}_\infty
 \quad\text{on }M_\infty.
\]
The metric $\mathbf{h}_\infty$ is complete by the static compactness
theorem; the last assertion of
Lemma~\ref{lem:per-distorsion-metrica-distancias} implies that
$\mathbf{g}_\infty(t)$ is also complete.
\end{proof}

\section{Rescalings and singularity models}
\label{sec:per-reescalamientos-singularidad}

The preceding compactness theorem becomes useful once the correct scale is chosen.
At a point where $|\operatorname{Rm}|=Q$, the metric $Q\mathbf{g}$ makes the
curvature have size one. Time must be multiplied by the same
factor to preserve the Ricci equation. The following formulas fix
all conventions used for singularity limits.

\begin{definition}[Temporal types of solutions]
\label{def:per-tipos-soluciones-ricci}
Let $M$ be a smooth manifold without boundary and let $\mathbf{g}(t)$ be a solution
of Ricci flow on $M$. For a problem with initial data at
$t=0$, the solution is \emph{forward maximal} if it admits no
smooth extension beyond its right time endpoint. This is the
meaning of maximality used for flows on $[0,T)$.
Independently of this property, a solution is called \emph{ancient} if its time interval
has left endpoint $-\infty$, for example $(-\infty,T]$;
\emph{eternal} if defined on $\mathbb R$, and \emph{immortal} if defined
on $[0,\infty)$. A \emph{rescaling limit model} is
a nonflat pointed smooth limit of parabolic rescalings about
points whose curvature tends to infinity.
\end{definition}

Suppose the original flow is defined on $[0,T)$. Let
$(x_i,t_i)\in M\times[0,T)$, with $i\in\mathbb N$, and set
\[
Q_i:=|\operatorname{Rm}_{\mathbf{g}(t_i)}|_{\mathbf{g}(t_i)}(x_i)>0.
\]
The parabolic rescaling centered at $(x_i,t_i)$ is
\begin{equation}
\mathbf{g}_i(s):=Q_i\,\mathbf{g}\!\left(t_i+\frac{s}{Q_i}\right),
\qquad
s\in\bigl[-Q_it_i,Q_i(T-t_i)\bigr).
\label{eq:per-reescalamiento-parabolico}
\end{equation}

For each $i\in\mathbb N$ and $q\in\mathbb N_0$, the scaling laws are
\begin{align}
s&=Q_i(t-t_i),&
\mathbf{g}_i&=Q_i \mathbf{g},&
d_{\mathbf{g}_i}&=Q_i^{1/2}d_{\mathbf{g}},\nonumber\\
d\lambda_{\mathbf{g}_i}&=Q_i^{m/2}d\lambda_{\mathbf{g}},&
\operatorname{inj}_{\mathbf{g}_i}&=Q_i^{1/2}\operatorname{inj}_{\mathbf{g}},&
R_{\mathbf{g}_i}&=Q_i^{-1}R_{\mathbf{g}},\nonumber\\
\operatorname{Rm}^{(1,3)}_{\mathbf{g}_i}
&=\operatorname{Rm}^{(1,3)}_{\mathbf{g}},&
\operatorname{Rm}^{(0,4)}_{\mathbf{g}_i}
&=Q_i\operatorname{Rm}^{(0,4)}_{\mathbf{g}},&
\operatorname{Ric}_{\mathbf{g}_i}&=\operatorname{Ric}_{\mathbf{g}},\nonumber\\
|\nabla_{\mathbf{g}_i(s)}^q\operatorname{Rm}_{\mathbf{g}_i(s)}|_{\mathbf{g}_i(s)}
&=Q_i^{-1-q/2}
|\nabla_{\mathbf{g}(t_i+s/Q_i)}^q
\operatorname{Rm}_{\mathbf{g}(t_i+s/Q_i)}|_{\mathbf{g}(t_i+s/Q_i)}.
\label{eq:per-leyes-reescalamiento}
\end{align}
In particular,
$|\operatorname{Rm}_{\mathbf{g}_i(0)}|_{\mathbf{g}_i(0)}(x_i)=1$. Moreover,
\[
\partial_s\mathbf{g}_i
=Q_i\frac1{Q_i}\partial_t\mathbf{g}\!\left(t_i+\frac{s}{Q_i}\right)
=-2\operatorname{Ric}_{\mathbf{g}(t_i+s/Q_i)}
=-2\operatorname{Ric}_{\mathbf{g}_i},
\]
so the flow is invariant.

The choice of points cannot be arbitrary. Take an increasing sequence
$(t_i)_{i\in\mathbb N}$ with $t_i<T$ for each $i\in\mathbb N$. If
\begin{equation}
Q_i=|\operatorname{Rm}_{\mathbf{g}(t_i)}|_{\mathbf{g}(t_i)}(x_i)
=\max_{(x,t)\in M\times[0,t_i]}|\operatorname{Rm}_{\mathbf{g}(t)}|_{\mathbf{g}(t)}(x),
\qquad t_i\to T^{-},
\label{eq:per-maximos-espaciotemporales}
\end{equation}
then
\[
|\operatorname{Rm}_{\mathbf{g}_i(s)}|_{\mathbf{g}_i(s)}\leq1
\quad\text{for }-Q_it_i\leq s\leq0,
\qquad
|\operatorname{Rm}_{\mathbf{g}_i(0)}|_{\mathbf{g}_i(0)}(x_i)=1.
\]
Shi estimates control all derivatives on compact subintervals
of $(-\infty,0]$. A lower injectivity radius bound is still
needed; we will obtain it from Perelman's noncollapsing. Hamilton
compactness will then produce a complete nonflat ancient solution.

\section{Heat with a time-dependent metric}
\label{sec:per-calor-conjugado}

Perelman's monotone quantities are constructed by coupling Ricci
flow with a backward heat equation. Since Riemannian
measure is also time-dependent, the spacetime adjoint is not
$-\partial_t-\Delta$; a scalar curvature term appears to compensate
for volume variation.

For a Ricci flow define
\[
\Box:=\partial_t-\Delta_{\mathbf{g}(t)}.
\]
The measure evolution, proved in
\eqref{eq:evolucion-medida-ricci}, is
$\partial_td\lambda_{\mathbf{g}(t)}=-R_{\mathbf{g}(t)}\,d\lambda_{\mathbf{g}(t)}$.

\begin{proposition}[Spacetime adjoint]
\label{prop:per-adjunto-calor-conjugado}
If $u,v$ are smooth and $M$ is closed, then
\begin{equation}
\int_{t_1}^{t_2}\!\!\int_M v\,\Box u\,d\lambda_{\mathbf{g}(t)}dt
=
\left[\int_Muv\,d\lambda_{\mathbf{g}(t)}\right]_{t_1}^{t_2}
+\int_{t_1}^{t_2}\!\!\int_Mu\,\Box^*v\,d\lambda_{\mathbf{g}(t)}dt,
\label{eq:per-integracion-calor-conjugado}
\end{equation}
where
\begin{equation}
\Box^*:=-\partial_t-\Delta_{\mathbf{g}(t)}+R_{\mathbf{g}(t)}.
\label{eq:per-operador-calor-conjugado}
\end{equation}
Thus the conjugate equation is
\[
\Box^*u=0,
\qquad
\partial_tu=-\Delta_{\mathbf{g}(t)}u+R_{\mathbf{g}(t)}u.
\]
Every positive solution of this equation preserves
$\displaystyle \int_Mu\,d\lambda_{\mathbf{g}(t)}$.
\end{proposition}

\begin{proof}
Differentiating the time pairing gives
\[
\frac d{dt}\int_Muv\,d\lambda_{\mathbf{g}(t)}
=\int_M(u_tv+uv_t-R_{\mathbf{g}(t)}uv)\,d\lambda_{\mathbf{g}(t)}.
\]
Spatial self-adjointness of $\Delta_{\mathbf{g}(t)}$ and integration in time give
\eqref{eq:per-integracion-calor-conjugado}. If $\Box^*u=0$, then
\[
\frac d{dt}\int_Mu\,d\lambda_{\mathbf{g}(t)}
=\int_M\bigl(-\Delta_{\mathbf{g}(t)}u+R_{\mathbf{g}(t)}u-R_{\mathbf{g}(t)}u\bigr)
\,d\lambda_{\mathbf{g}(t)}=0.
\]
\end{proof}

Time dependence prevents use of an autonomous semigroup. For
$s\leq t$ write $U(t,s)$ for the propagator of
$\partial_tu=\Delta_{\mathbf{g}(t)}u$. It satisfies
\[
U(t,r)U(r,s)=U(t,s),
\qquad s\leq r\leq t,
\]
but does not depend only on $t-s$. Its two-time kernel
$H(x,t;y,s)$ satisfies the same composition law, integrating over the
intermediate slice. The adjoint of $U(t,s)$ propagates backward through
$\Box^*$; in general there is no symmetry under simultaneous interchange of
$(x,t)$ and $(y,s)$. The construction follows by applying the parabolic
theory of Chapter~\ref{cap:teoria-parabolica-cerradas} on short
intervals and composing the evolution operators.

\section{The functional \texorpdfstring{$\mathcal F$}{F}}
\label{sec:per-funcional-F}

Perelman's variational description considers a metric and a density
simultaneously. The auxiliary function $f$ allows the weighted measure to be fixed
and introduces an energy whose derivative along the flow is a square.

Let $M$ be closed. For $f\in C^\infty(M)$ define
\begin{equation}
 \mathcal F(\mathbf{g},f)
 :=\int_M\bigl(R_{\mathbf{g}}+|\nabla f|_{\mathbf{g}}^2\bigr)e^{-f}
 \,d\lambda_{\mathbf{g}}.
 \label{eq:per-definicion-F}
\end{equation}
The natural normalization is
$\displaystyle \int_Me^{-f}\,d\lambda_{\mathbf{g}}=1$.

Before differentiating $\mathcal F$, collect the integrations by parts
explaining the Hessian terms that will appear.

\begin{lemma}[Integration by parts with the weighted measure]
\label{lem:per-integracion-ponderada-F}
Let $M$ be a closed smooth manifold, let $f,\phi\in C^\infty(M)$, and let
$\mathbf{h}$ be a symmetric covariant tensor of order two. Then
\begin{align}
 \int_M\operatorname{div}_{\mathbf{g}}^2\mathbf{h}\,e^{-f}
 \,d\lambda_{\mathbf{g}}
 &=\int_M\bigl[-\langle\mathbf{h},\nabla_{\mathbf{g}}^2f\rangle_{\mathbf{g}}
 +\mathbf{h}(\nabla^{\mathbf{g}}f,\nabla^{\mathbf{g}}f)\bigr]e^{-f}
 \,d\lambda_{\mathbf{g}},
 \label{eq:per-partes-div2-ponderada}\\
 -\int_M\Delta_{\mathbf{g}}(\operatorname{tr}_{\mathbf{g}}\mathbf{h})e^{-f}
 \,d\lambda_{\mathbf{g}}
 &=\int_M\operatorname{tr}_{\mathbf{g}}\mathbf{h}
 \bigl(\Delta_{\mathbf{g}}f-|\nabla f|_{\mathbf{g}}^2\bigr)e^{-f}
 \,d\lambda_{\mathbf{g}},
 \label{eq:per-partes-traza-ponderada}\\
 2\int_M\langle\nabla^{\mathbf{g}}f,\nabla^{\mathbf{g}}\phi\rangle_{\mathbf{g}}
 e^{-f}\,d\lambda_{\mathbf{g}}
 &=-2\int_M\phi\bigl(\Delta_{\mathbf{g}}f-|\nabla f|_{\mathbf{g}}^2\bigr)e^{-f}
 \,d\lambda_{\mathbf{g}}.
 \label{eq:per-partes-phi-ponderada}
\end{align}
\end{lemma}

\begin{proof}
For the first identity, integrate the two divergences successively:
\begin{align*}
 \int_M\sum_{i,j=1}^{m}\nabla^i\nabla^jh_{ij}\,e^{-f}\,d\lambda_{\mathbf{g}}
 &=\int_M\sum_{i,j=1}^{m}(\nabla^jh_{ij})\nabla^if\,e^{-f}\,d\lambda_{\mathbf{g}}\\
 &=-\int_M\sum_{i,j=1}^{m}h_{ij}\nabla^j(\nabla^if\,e^{-f})\,d\lambda_{\mathbf{g}}.
\end{align*}
Since
$\nabla^j(\nabla^if\,e^{-f})=(\nabla^2f)^{ij}e^{-f}
-\nabla^if\nabla^jf\,e^{-f}$, we obtain
\eqref{eq:per-partes-div2-ponderada}. If
$q=\operatorname{tr}_{\mathbf{g}}\mathbf{h}$, then
\[
 -\int_M(\Delta_{\mathbf{g}}q)e^{-f}\,d\lambda_{\mathbf{g}}
 =-\int_M\langle\nabla q,\nabla f\rangle_{\mathbf{g}}e^{-f}
 \,d\lambda_{\mathbf{g}}.
\]
Integrating the divergence of $q e^{-f}\nabla^{\mathbf{g}}f$ gives
\eqref{eq:per-partes-traza-ponderada}. The same calculation with $q=2\phi$
proves \eqref{eq:per-partes-phi-ponderada}.
\end{proof}

\begin{proposition}[First variation of $\mathcal F$]
\label{prop:per-primera-variacion-F}
Let $M$ be a closed smooth manifold. For a metric $\mathbf{g}$, a
function $f\in C^\infty(M)$, and variations
$\mathbf{h}=\delta\mathbf{g}$, $\phi=\delta f$, we have
\begin{align}
 D\mathcal F_{(\mathbf{g},f)}[\mathbf{h},\phi]
 =\int_M\Bigg[-\langle\mathbf{h},
 \operatorname{Ric}_{\mathbf{g}}+\nabla_{\mathbf{g}}^2f\rangle_{\mathbf{g}}+\left(\frac12\operatorname{tr}_{\mathbf{g}}\mathbf{h}-\phi\right)
 \bigl(2\Delta_{\mathbf{g}}f-|\nabla f|_{\mathbf{g}}^2+R_{\mathbf{g}}\bigr)
 \Bigg]e^{-f}\,d\lambda_{\mathbf{g}}.
 \label{eq:per-primera-variacion-F}
\end{align}
\end{proposition}

\begin{proof}
The variation formulas of the preceding chapter give
\begin{align*}
 \delta R_{\mathbf{g}}
 &=-\langle\mathbf{h},\operatorname{Ric}_{\mathbf{g}}\rangle_{\mathbf{g}}
 +\operatorname{div}_{\mathbf{g}}^2\mathbf{h}
 -\Delta_{\mathbf{g}}\operatorname{tr}_{\mathbf{g}}\mathbf{h},\\
 \delta|\nabla f|_{\mathbf{g}}^2
 &=-\mathbf{h}(\nabla^{\mathbf{g}}f,\nabla^{\mathbf{g}}f)
 +2\langle\nabla^{\mathbf{g}}f,\nabla^{\mathbf{g}}\phi\rangle_{\mathbf{g}},\\
 \delta(e^{-f}d\lambda_{\mathbf{g}})
 &=\left(\frac12\operatorname{tr}_{\mathbf{g}}\mathbf{h}-\phi\right)
 e^{-f}d\lambda_{\mathbf{g}}.
\end{align*}
Substituting into \eqref{eq:per-definicion-F}, the sum of terms
containing derivatives of $\mathbf{h}$ becomes, by
\eqref{eq:per-partes-div2-ponderada} and
\eqref{eq:per-partes-traza-ponderada},
\[
 -\langle\mathbf{h},\nabla^2f\rangle
 +\mathbf{h}(\nabla f,\nabla f)
 +\operatorname{tr}_{\mathbf{g}}\mathbf{h}
   (\Delta f-|\nabla f|^2).
\]
The second summand cancels the metric variation in
$|\nabla f|^2$ exactly. By \eqref{eq:per-partes-phi-ponderada}, the term
containing $\nabla\phi$ becomes
$-2\phi(\Delta f-|\nabla f|^2)$. Adding the measure variation, the
coefficients of $\operatorname{tr}_{\mathbf{g}}\mathbf{h}$ and $\phi$ are,
respectively,
\[
 \frac12(2\Delta f-|\nabla f|^2+R),
 \qquad -(2\Delta f-|\nabla f|^2+R).
\]
This proves \eqref{eq:per-primera-variacion-F}.
\end{proof}

\begin{theorem}[Monotonicity of $\mathcal F$]
\label{teo:per-monotonia-F}
Let $M$ be a closed smooth manifold and suppose that
\begin{equation}
 \partial_t\mathbf{g}=-2\operatorname{Ric}_{\mathbf{g}(t)},
 \qquad
 \partial_tf=-\Delta_{\mathbf{g}(t)}f+|\nabla f|_{\mathbf{g}(t)}^2
 -R_{\mathbf{g}(t)}.
 \label{eq:per-sistema-F}
\end{equation}
Then $u=e^{-f}$ satisfies $\Box^*u=0$ and
\begin{equation}
 \frac d{dt}\mathcal F(\mathbf{g}(t),f(t))
 =2\int_M|\operatorname{Ric}_{\mathbf{g}(t)}+
 \nabla_{\mathbf{g}(t)}^2f|_{\mathbf{g}(t)}^2e^{-f}
 \,d\lambda_{\mathbf{g}(t)}\geq0.
 \label{eq:per-monotonia-F}
\end{equation}
Equality on an interval characterizes, modulo diffeomorphisms, a steady
gradient soliton. On a closed manifold, such a soliton is
Ricci-flat and $f$ is constant on every connected component.
\end{theorem}

\begin{proof}
The equation for $f$ implies
\[
 \partial_t(e^{-f})
 =(\Delta f-|\nabla f|^2+R)e^{-f}
 =-\Delta(e^{-f})+Re^{-f},
\]
so $\Box^*u=0$. Let $\varphi_t$ be the diffeomorphisms generated by
$-\nabla^{\mathbf{g}(t)}f$. By naturality of Ricci and the Hessian, the
pullback of $(\mathbf{g},f)$ satisfies
\[
 \partial_t\mathbf{g}=-2(\operatorname{Ric}_{\mathbf{g}}+\nabla_{\mathbf{g}}^2f),
 \qquad \partial_tf=-\Delta_{\mathbf{g}}f-R_{\mathbf{g}}.
\]
In this gauge,
$\frac12\operatorname{tr}_{\mathbf{g}}(\partial_t\mathbf{g})-\partial_tf=0$;
thus the measure $e^{-f}d\lambda_{\mathbf{g}}$ remains fixed
pointwise. Apply \eqref{eq:per-primera-variacion-F} with
$\mathbf{h}=-2(\operatorname{Ric}+\nabla^2f)$ and
$\phi=-\Delta f-R$. The second variation term vanishes, and the first
is the right-hand side of \eqref{eq:per-monotonia-F}. Diffeomorphism invariance of
$\mathcal F$ returns the formula to the original flow.

If equality holds, then
$\operatorname{Ric}+\nabla^2f=0$. Taking the trace and integrating gives
$\displaystyle \int_MR\,d\lambda_{\mathbf{g}}=0$. Naturality of scalar curvature
and its evolution equation give
\[
 \Delta_fR:=\Delta R-\langle\nabla f,\nabla R\rangle
 =-2|\operatorname{Ric}|^2.
\]
The measure $e^{-f}d\lambda_{\mathbf{g}}$ makes $\Delta_f$ symmetric, so
the integral of the left-hand side is zero. We conclude that
$\operatorname{Ric}=0$ and then $\nabla^2f=0$. A function with zero Hessian on a closed manifold is constant on every
connected component.
\end{proof}

\section{The entropy \texorpdfstring{$\mathcal W$}{W} and the quantity
\texorpdfstring{$\mu$}{mu}}
\label{sec:per-W-mu}

The functional $\mathcal F$ is not invariant under parabolic
rescaling. The entropy $\mathcal W$ incorporates a time scale
$\tau>0$ and a Boltzmann entropy term so that the entire
expression remains invariant when $\mathbf{g}$ and $\tau$ are multiplied
by the same constant.

Let
\[
 u=(4\pi\tau)^{-m/2}e^{-f},
 \qquad \int_Mu\,d\lambda_{\mathbf{g}}=1,
\]
and define
\begin{equation}
 \mathcal W(\mathbf{g},f,\tau)
 :=\int_M\left[\tau(R_{\mathbf{g}}+|\nabla f|_{\mathbf{g}}^2)+f-m\right]
 u\,d\lambda_{\mathbf{g}}.
 \label{eq:per-definicion-W}
\end{equation}

\begin{lemma}[Boltzmann entropy and $\mathcal F$]
\label{lem:per-boltzmann-F}
Suppose that $\mathbf{g}(t)$ solves Ricci flow and a positive function
$u$ satisfies $\Box^*u=0$ and
$\displaystyle \int_Mu\,d\lambda_{\mathbf{g}(t)}=1$. If
\[
 \mathcal B(t):=\int_Mu\log u\,d\lambda_{\mathbf{g}(t)},
 \qquad
 \mathcal F(t):=\int_M\left(\frac{|\nabla u|^2}{u}+Ru\right)
 d\lambda_{\mathbf{g}(t)},
\]
then $\mathcal B'(t)=\mathcal F(t)$.
\end{lemma}

\begin{proof}
Use $u_t=-\Delta u+Ru$ and
$\partial_td\lambda_{\mathbf{g}}=-R\,d\lambda_{\mathbf{g}}$. Differentiation under
the integral gives
\begin{align*}
 \mathcal B'(t)
 =\int_M\bigl[u_t(1+\log u)-Ru\log u\bigr]d\lambda_{\mathbf{g}}=\int_M\bigl[-\Delta u(1+\log u)+Ru\bigr]d\lambda_{\mathbf{g}}.
\end{align*}
The integral of $\Delta u$ is zero and
$\displaystyle -\int_M(\Delta u)\log u=\int_M|\nabla u|^2/u$. This proves the identity.
\end{proof}

\begin{theorem}[Monotonicity of $\mathcal W$]
\label{teo:per-monotonia-W}
Let $M$ be a closed smooth manifold of dimension $m$. Suppose that
\begin{equation}
 \partial_t\mathbf{g}=-2\operatorname{Ric}_{\mathbf{g}(t)},
 \qquad \tau'=-1,
 \qquad
 \partial_tf=-\Delta_{\mathbf{g}(t)}f+|\nabla f|_{\mathbf{g}(t)}^2
 -R_{\mathbf{g}(t)}+\frac m{2\tau},
 \label{eq:per-sistema-W}
\end{equation}
and $\tau(t)>0$. Also suppose that the density
$u=(4\pi\tau)^{-m/2}e^{-f}$ has mass one at some time in the interval.
Then $u$ solves $\Box^*u=0$, preserves this normalization, and
\begin{equation}
 \frac d{dt}\mathcal W(\mathbf{g}(t),f(t),\tau(t))
 =2\tau\int_M\left|
 \operatorname{Ric}_{\mathbf{g}(t)}+\nabla_{\mathbf{g}(t)}^2f
 -\frac1{2\tau}\mathbf{g}(t)
 \right|_{\mathbf{g}(t)}^2u\,d\lambda_{\mathbf{g}(t)}\geq0.
 \label{eq:per-monotonia-W}
\end{equation}
Equality on an interval is equivalent to the equation of a shrinking
gradient soliton.
\end{theorem}

\begin{proof}
From
$\log u=-\frac m2\log(4\pi\tau)-f$ and \eqref{eq:per-sistema-W} we obtain
$u_t=-\Delta u+Ru$. Mass is preserved by
Proposition~\ref{prop:per-adjunto-calor-conjugado}.

Write
\[
 \mathcal F_u(t):=\int_M(R+|\nabla f|^2)u\,d\lambda_{\mathbf{g}},
 \qquad
 \mathcal B(t):=\int_Mu\log u\,d\lambda_{\mathbf{g}}.
\]
Since $u=(4\pi\tau)^{-m/2}e^{-f}$ and its mass is one,
\begin{equation}
 \mathcal W=\tau\mathcal F_u-\mathcal B
 -\frac m2\log(4\pi\tau)-m.
 \label{eq:per-W-F-B}
\end{equation}
Lemma~\ref{lem:per-boltzmann-F} gives $\mathcal B'=\mathcal F_u$. To
calculate the derivative of the first term, we do not apply
Theorem~\ref{teo:per-monotonia-F} directly to $f$, since its equation contains the
summand $m/(2\tau)$. Define
\[
 \widetilde f:=f+\frac m2\log(4\pi\tau)=-\log u.
\]
Since $\tau'=-1$, we have
\[
 \partial_t\widetilde f
 =-\Delta\widetilde f+|\nabla\widetilde f|^2-R,
 \qquad e^{-\widetilde f}=u.
\]
Moreover, $d\widetilde f=df$ and $\nabla^2\widetilde f=\nabla^2f$. Thus,
\[
 \mathcal F_u=\mathcal F(\mathbf{g},\widetilde f),
 \qquad
 \mathcal F_u'=2\int_M|\operatorname{Ric}+\nabla^2f|^2u
 \,d\lambda_{\mathbf{g}}
\]
by Theorem~\ref{teo:per-monotonia-F}.
Differentiating \eqref{eq:per-W-F-B} and using $\tau'=-1$ gives
\begin{equation}
 \mathcal W'
 =2\tau\int_M|\operatorname{Ric}+\nabla^2f|^2u\,d\lambda_{\mathbf{g}}
 -2\mathcal F_u+\frac m{2\tau}.
 \label{eq:per-W-antes-cuadrado}
\end{equation}
Weighted integration by parts gives
\[
 \int_M(\Delta f)u\,d\lambda_{\mathbf{g}}
 =\int_M|\nabla f|^2u\,d\lambda_{\mathbf{g}};
\]
consequently,
$\displaystyle \mathcal F_u=\int_M(R+\Delta f)u
=\int_M\operatorname{tr}_{\mathbf{g}}(\operatorname{Ric}+\nabla^2f)u$.
If $\mathbf{A}=\operatorname{Ric}+\nabla^2f$, then
\[
 2\tau\left|\mathbf{A}-\frac1{2\tau}\mathbf{g}\right|^2
 =2\tau|\mathbf{A}|^2-2\operatorname{tr}_{\mathbf{g}}\mathbf{A}
 +\frac m{2\tau}.
\]
Integrating this identity and comparing it with
\eqref{eq:per-W-antes-cuadrado} proves
\eqref{eq:per-monotonia-W}. Equality is equivalent to
$\operatorname{Ric}+\nabla^2f=\mathbf{g}/(2\tau)$.
\end{proof}

\begin{definition}[The quantity $\mu$]
\label{def:per-mu}
Let $M$ be a closed smooth manifold of dimension $m$, let $\mathbf{g}$ be a Riemannian
metric, and let $\tau>0$. Set
\begin{equation}
 \mu(\mathbf{g},\tau)
 :=\inf\left\{\mathcal W(\mathbf{g},f,\tau)\middle|
 f\in C^\infty(M),\quad
 (4\pi\tau)^{-m/2}\int_Me^{-f}\,d\lambda_{\mathbf{g}}=1\right\}.
 \label{eq:per-definicion-mu}
\end{equation}
\end{definition}

\begin{lemma}[Strong principle for the entropy nonlinearity]
\label{lem:per-maximo-entropico}
Let $\Omega$ be a connected domain in a Riemannian manifold without boundary and
let $w\in C^2(\Omega)$ be nonnegative. Suppose that there exists $s_0>0$ such that,
when $0\leq w\leq s_0$,
\[
 \Delta w\leq\beta(w),
 \qquad \beta(s)=C_0s-C_1s\log s^2,
 \quad C_1>0.
\]
Then $w$ is identically zero or strictly positive.
The same alternative holds for a nonnegative function
$w\in W^{2,p_0}_{\mathrm{loc}}(\Omega)$ with
$p_0>\displaystyle\max\{2,\dim\Omega\}$, taking its $C^1$ representative, if the
inequality holds almost everywhere on $\{0\leq w\leq s_0\}$.
\end{lemma}

\begin{proof}
First prove the assertion for $w\in C^2(\Omega)$.
Extend $\beta$ by $\beta(0)=0$ and reduce $s_0>0$ if
necessary. Then $\beta$ is continuous, nonnegative, and increasing on
$[0,s_0]$, because
$\beta'(s)=C_0-C_1(\log s^2+2)>0$ for sufficiently small $s$.
We may also assume that $s_0<1$ and $\beta$ is strictly increasing.

Suppose $w$ has both a zero and a positive point. Connectedness gives
$z_0\in\partial\{w>0\}\cap\Omega$; continuity implies $w(z_0)=0$.
Let $d$ be the Riemannian distance on $\Omega$.
Lemma~\ref{lem:vecindades-normales-convexas-locales}, applied at $z_0$ to
an open set $O\Subset\Omega$, gives an open set $W\subseteq O$ and
$\varepsilon_0>0$ such that, for each $p\in W$, $\exp_p$ is a
diffeomorphism on $B_{\mathbf g_p}(0,2\varepsilon_0)$ and
$W\subseteq\exp_p(B_{\mathbf g_p}(0,\varepsilon_0/2))$.
Choose an open set $V$ with $z_0\in V$ and compact closure contained
in $W$, then choose $0<a<\varepsilon_0/4$ with $B_d(z_0,4a)\subset V$.
The closed ball $\overline B_d(z_0,3a)$ is a closed subset of
$\overline V$ and is therefore compact.

Take $p\in\{w>0\}\cap B_d(z_0,a)$ and set
$Z:=\{x\in\Omega\mid w(x)=0\}$ and $R:=d(p,Z)$. Continuity of $w$
implies $R>0$, while $z_0\in Z$ gives $R\leq d(p,z_0)<a$.
If $x\in Z\setminus\overline B_d(z_0,3a)$, then
\[
 d(p,x)\geq d(z_0,x)-d(z_0,p)>2a>R.
\]
Thus the infimum defining $R$ is attained on the nonempty compact set
$Z\cap\overline B_d(z_0,3a)$. There is a point $q$ in that set with
$d(p,q)=R$. Moreover,
\[
 \overline B_d(p,R)\subset B_d(z_0,2a)\subset V\subset W,
 \qquad w>0\quad\text{on }B_d(p,R),\qquad w(q)=0.
\]
The closed ball centered at $p$ is compact and lies in the chosen normal
domain. By radial minimization in
Proposition~\ref{prop:minimizacion-radial-vecindad-normal}, throughout this
ball $d(p,x)=|\exp_p^{-1}(x)|_{\mathbf g_p}$. Hence
the function $r(x)=d(p,x)$ is smooth on the annulus
$\mathcal A=\{R/2<r<R\}$. Fix $C\geq0$ such that $\Delta r\leq C$
on its closure and set $L=R/2$. For $0<\epsilon<1$, define
\[
 k_\epsilon=C+
 \sqrt{|C_0|+2C_1|\log\epsilon|+2\frac{C_1}{e}+1},
 \qquad
 z_\epsilon(t)=\epsilon(e^{k_\epsilon t}-1),
 \quad 0\leq t\leq L.
\]
Since $k_\epsilon=O(\sqrt{|\log\epsilon|})$,
$z_\epsilon(L)\to0$ as $\epsilon\to0^{+}$.
Choose $\epsilon$ small enough that
\[
 0<z_\epsilon(L)<
 \min\left\{s_0,\min_{r=R/2}w\right\}.
\]
Let $\phi(x)=z_\epsilon(R-r(x))$. On both components of
$\partial\mathcal A$ we have $\phi\leq w$.

Let us verify the differential inequality for the barrier. For $t>0$,
write $y=e^{k_\epsilon t}-1>0$. Since
$\displaystyle\frac{y}{1+y}\displaystyle\max\{-\log y,0\}\leq1/e$, we obtain
\[
 \beta(z_\epsilon(t))
 \leq\epsilon e^{k_\epsilon t}
       \bigl(|C_0|+2C_1|\log\epsilon|+2C_1/e\bigr).
\]
On the other hand, $|dr|=1$ and the choice of $k_\epsilon$ give
\[
 \Delta\phi
 =z_\epsilon''(R-r)-z_\epsilon'(R-r)\Delta r
 \geq\epsilon e^{k_\epsilon(R-r)}
          (k_\epsilon^2-Ck_\epsilon)
 \geq\beta(\phi).
\]
If $\phi-w$ had a positive maximum in $\mathcal A$, at that point
we would have $0\leq w<\phi<s_0$ and
\[
 0\geq\Delta(\phi-w)\geq\beta(\phi)-\beta(w)>0,
\]
a contradiction. Thus $w\geq\phi$ in the annulus.
Let $\gamma(t)$ be the radial geodesic from $q$ toward $p$,
parametrized by arc length. For $0<t<L$ we have
$w(\gamma(t))\geq z_\epsilon(t)$, and both functions vanish at $t=0$.
Their right derivatives satisfy
\[
 dw_q(\dot\gamma(0))\geq z_\epsilon'(0)=\epsilon k_\epsilon>0.
\]
This contradicts $dw_q=0$, since $w\geq0$ and $w(q)=0$ at an
interior point of $\Omega$. This proves the alternative.

Now consider $w\in W^{2,p_0}_{\mathrm{loc}}(\Omega)$ with
$p_0>\displaystyle\max\{2,\dim\Omega\}$. The local Sobolev--Morrey embedding, obtained
by applying Theorem~\ref{teo:sobolev-morrey-euclidiano} in charts after
multiplication by cutoff functions, provides a representative of class $C^1$.
Suppose again that $w$ has a zero and a positive point.
The choice of the ball, the annulus, and the barrier $\phi$ uses only
the continuity of $w$ and remains valid. To compare $w$ with $\phi$,
set $h:=\phi-w$ and, given $\delta>0$,
\[
 \psi_\delta:=(h-\delta)_+.
\]
Since $h\leq0$ on $\partial\mathcal A$, the support of $\psi_\delta$ is
contained in the compact set
$\{x\in\overline{\mathcal A}\mid h(x)\geq\delta\}\subseteq\mathcal A$.
The chain rule in Theorem~\ref{teo:regla-cadena-lipschitz-sobolev} gives
$\psi_\delta\in W^{1,2}(\mathcal A)$ and
\[
 d\psi_\delta=\mathbf 1_{\{h>\delta\}}\,dh
 \qquad\text{almost everywhere.}
\]
We can approximate $\psi_\delta$ in $W^{1,2}(M)$ by smooth functions
with compact support in $\mathcal A$: apply
Theorem~\ref{meyers-serrin-haz} to the trivial bundle and multiply the
approximations by a fixed cutoff function with compact support in $\mathcal A$
that equals one near the support of $\psi_\delta$.
Moreover, $\Delta h\in L^{p_0}(\mathcal A)\subseteq L^2(\mathcal A)$,
because $p_0>2$ and the annulus has finite volume. Thus the weak
integration by parts identity admits $\psi_\delta$ as a test function.
Where $\psi_\delta>0$, we have $0\leq w<\phi<s_0$, and we obtain
\[
 -\int_{\mathcal A}|d\psi_\delta|_{\mathbf g}^2\,d\lambda_{\mathbf g}
 =\int_{\mathcal A}(\Delta h)\psi_\delta\,d\lambda_{\mathbf g}
 \geq\int_{\mathcal A}\bigl(\beta(\phi)-\beta(w)\bigr)
             \psi_\delta\,d\lambda_{\mathbf g}\geq0.
\]
The last integrand is strictly positive wherever $\psi_\delta>0$.
Thus this set has measure zero and, being open, is empty.
It follows that $h\leq\delta$ for every $\delta>0$, and hence $w\geq\phi$ in
$\mathcal A$. The contradiction with $dw_q=0$ now follows along
the same radial geodesic, since the representative of $w$ is of class $C^1$.
This proves the weak assertion.

For this nonlinearity, the barrier above implements the criterion in
Vázquez's strong principle~\cite{Vazquez1984}. Indeed, with
$\displaystyle B(s)=\int_0^s\beta(r)\,dr$,
\[
 \int_0^{s_0}\frac{ds}{\sqrt{B(s)}}=\infty.
\]
In our case, integration gives
$B(s)=(C_0+C_1)\frac{s^2}{2}-C_1s^2\log s$. Consequently,
$B(s)\leq Cs^2(1+|\log s|)$ near zero. Under the substitution
$z=-\log s$, the integral providing the lower bound is a multiple of
$\displaystyle \int_{-\log s_0}^{\infty}(1+z)^{-1/2}\,dz$, which diverges.
This divergence explains why the logarithmic term still permits
the strong principle, even though $\beta(s)/s$ is unbounded near zero.
\end{proof}

\begin{proposition}[Minimizers and invariances of $\mu$]
\label{prop:per-minimizador-mu}
On a smooth closed connected manifold, the infimum of
\eqref{eq:per-definicion-mu} is attained by a smooth function. If
\[
 v=(4\pi\tau)^{-m/4}e^{-f/2},
 \qquad \int_Mv^2\,d\lambda_{\mathbf{g}}=1,
\]
the minimizer can be chosen positive and satisfies
\begin{equation}
 -4\tau\Delta_{\mathbf{g}}v+\tau R_{\mathbf{g}}v-v\log v^2
 -\frac m2\log(4\pi\tau)v-mv=\mu(\mathbf{g},\tau)v.
 \label{eq:per-euler-mu-v}
\end{equation}
Equivalently,
\begin{equation}
 2\tau\Delta_{\mathbf{g}}f-\tau|\nabla f|_{\mathbf{g}}^2
 +\tau R_{\mathbf{g}}+f-m=\mu(\mathbf{g},\tau).
 \label{eq:per-euler-mu-f}
\end{equation}
Moreover, for every diffeomorphism $\varphi$ and every $c>0$,
\begin{equation}
 \mu(\varphi^*\mathbf{g},\tau)=\mu(\mathbf{g},\tau),
 \qquad \mu(c\mathbf{g},c\tau)=\mu(\mathbf{g},\tau).
 \label{eq:per-invariancias-mu}
\end{equation}
\end{proposition}

\begin{proof}
In the variable $v$, the problem is to minimize
\begin{equation}
 \mathcal J_{\mathbf{g},\tau}(v)
 :=\int_M\left[4\tau|\nabla v|_{\mathbf{g}}^2+\tau R_{\mathbf{g}}v^2
 -v^2\log v^2-\frac m2\log(4\pi\tau)v^2-mv^2\right]
 d\lambda_{\mathbf{g}}
 \label{eq:per-funcional-J-mu}
\end{equation}
over the sphere in $H^1(M)$ given by $\|v\|_{L^2(M)}=1$, with the
convention $0\log0=0$. Enlarging the admissible set in this way does not change
the infimum. To see this, approximate $|v|$ in $H^1(M)$ by nonnegative smooth
functions using the heat flow for the fixed metric; add a positive
constant tending to zero and normalize in $L^2(M)$. This gives
strictly positive smooth functions converging to $|v|$ in $H^1(M)$.
By the Sobolev embeddings, they also converge in some $L^q(M)$ with $q>2$.
The power estimate detailed below and Vitali's theorem
give convergence of their entropy terms. The quadratic
terms converge directly. Finally,
$\mathcal J(|v|)=\mathcal J(v)$, since
$|d|v||=|dv|$ almost everywhere for real functions in $H^1(M)$.

Fix $q>2$ for which $H^1(M)\hookrightarrow L^q(M)$ is continuous.
Jensen's inequality, applied to the probability measure
$v^2d\lambda_{\mathbf{g}}$, gives
\[
 \int_Mv^2\log v^2\,d\lambda_{\mathbf{g}}
 \leq\frac{2q}{q-2}\log\|v\|_{L^q(M)}.
\]
The Sobolev embedding and the inequality
$\log(1+a)\leq\varepsilon a^2+C_\varepsilon$ imply, for every
$\varepsilon>0$,
\begin{equation}
 \int_Mv^2\log v^2\,d\lambda_{\mathbf{g}}
 \leq\varepsilon\|dv\|_{L^2(M,T^*M)}^2+C_\varepsilon
 \qquad(\|v\|_{L^2(M)}=1).
 \label{eq:per-log-sobolev-coerciva-mu}
\end{equation}
Taking $0<\varepsilon<4\tau$ and using the boundedness of $R_{\mathbf{g}}$,
we obtain a lower bound for $\mathcal J_{\mathbf{g},\tau}$ and a
uniform bound in $H^1(M)$ for every minimizing sequence $(v_j)$.

We may replace $v_j$ by $|v_j|$. After passing to a subsequence,
\[
 v_j\rightharpoonup v\quad\text{in }H^1(M),
 \qquad v_j\to v\quad\text{in }L^{q_0}(M)
\]
for some $2<q_0<q$, and also almost everywhere. In particular,
$\|v\|_{L^2(M)}=1$. The Dirichlet energy is lower semicontinuous, and
the curvature term converges by strong convergence in $L^2(M)$.

For the entropy term, choose $\eta\in(0,1)$ with
$2(1+\eta)<q_0$. The inequality
\[
 a|\log a|\leq C_\eta(a^{1-\eta}+a^{1+\eta}),
 \qquad a\geq0,
\]
applied to $a=v_j^2$ shows that $(v_j^2\log v_j^2)$ is uniformly
integrable: the higher power is bounded in a space $L^p(M)$ with $p>1$ by
convergence in $L^{q_0}(M)$, and the lower power is bounded in
$L^{1/(1-\eta)}(M)$ by normalization in $L^2(M)$. Vitali's theorem gives
\[
 \int_Mv_j^2\log v_j^2\,d\lambda_{\mathbf{g}}
 \longrightarrow\int_Mv^2\log v^2\,d\lambda_{\mathbf{g}}.
\]
Therefore, $v$ is a nonnegative minimizer.

The real function $s\mapsto s^2\log s^2$, extended by zero at $s=0$, is
$C^1$, and its derivative grows more slowly than any power
$|s|^{1+\varepsilon}$ at infinity. The preceding embeddings allow us to
differentiate \eqref{eq:per-funcional-J-mu} in directions in $H^1(M)$. The
method of Lagrange multipliers yields
\[
 -4\tau\Delta v+\tau Rv-v\log v^2-v
 -\frac m2\log(4\pi\tau)v-mv=\lambda v
\]
in the weak sense. Multiplying by $v$ and integrating gives
$\lambda=\mu(\mathbf{g},\tau)-1$, and hence
\eqref{eq:per-euler-mu-v}.

Let us establish the regularity before using positivity. The weak equation
has the form
\[
 -\Delta v=F(x,v),\qquad
 |F(x,a)|\leq C_{\varepsilon}(1+|a|^{1+\varepsilon}),
 \quad x\in M,\quad a\in\mathbb R,
\]
for each $\varepsilon>0$. The term $a\log a^2$, extended by zero,
is bounded near $a=0$, which explains the constant summand.

If $m\geq3$, we start with
$v\in L^{q_0}(M)$ with $q_0=2m/(m-2)$. Fix
$\displaystyle 0<\varepsilon<\min\{q_0-1,q_0/m\}$. If $v\in L^{q_j}(M)$,
then $F(\cdot,v)\in L^{p_j}(M)$, where $p_j=q_j/(1+\varepsilon)>1$,
and elliptic regularity gives $v\in W^{2,p_j}(M)$.
As long as $p_j<m/2$, the Sobolev embedding yields the next exponent:
\[
 \frac1{q_{j+1}}=\frac{1+\varepsilon}{q_j}-\frac2m,
 \qquad j\in\mathbb N_0\text{ with }p_j<m/2.
\]
Since $q_j\geq q_0$, the choice of $\varepsilon$ implies
$\varepsilon/q_j<1/m$ and therefore
$1/q_{j+1}<1/q_j-1/m$. The reciprocals decrease by at least $1/m$ at
each step, so this alternative can occur only finitely
many times. If $p_j=m/2$ is reached, the critical embedding gives
$v\in L^q(M)$ for every $q<\infty$; we then choose
$q>(1+\varepsilon)m/2$ and apply the equation again. In either case,
we reach $W^{2,p}(M)$ with $p>m/2$, and hence $C^{0,\alpha}$ for some
$\alpha\in(0,1)$.

If $m=2$, the embedding of $H^1(M)$ gives $v\in L^q(M)$ for each
$q<\infty$. Choose $q>2(1+\varepsilon)$; the equation gives
$v\in W^{2,q/(1+\varepsilon)}(M)\hookrightarrow C^{1,\alpha}$ for some
$\alpha>0$. If $m=1$, the embedding of $H^1(M)$ already gives a continuous
bounded function. Thus $v$ is bounded in every dimension.

The function $F(\cdot,v)$ now belongs to $L^p(M)$ for each
$1<p<\infty$. Another application of elliptic regularity gives
$v\in W^{2,p}(M)$ for all these exponents and, taking $p>m$,
$v\in C^{1,\alpha}$. In particular, $v$ is Lipschitz.
Take $p>\displaystyle\max\{2,m\}$ to apply the weak version of the strong
principle just proved.

From \eqref{eq:per-euler-mu-v} we obtain, almost everywhere,
$\Delta v\leq C_0v-C_1v\log v^2$, where we may take
\[
 C_1=\frac1{4\tau},\qquad
 C_0=\frac1{4\tau}
 \left\|\tau R_{\mathbf{g}}-\frac m2\log(4\pi\tau)-m
                     -\mu(\mathbf{g},\tau)\right\|_\infty.
\]
The weak version of
Lemma~\ref{lem:per-maximo-entropico} and the normalization
exclude the zero solution; since $M$ is connected, it follows that $v>0$.
Compactness of $M$ gives $\displaystyle \min_{x\in M}v(x)>0$.
Thus $F$ is smooth on a neighborhood of the compact set
$\{(x,v(x))\mid x\in M\}$. We start with $v\in H^2(M)$.
If $v\in H^r(M)$ for an integer $r\geq2$, the composition estimate
\eqref{eq:moser-composicion-lineal-orden-alto} implies
$F(\cdot,v)\in H^r(M)$; the equation and
Theorem~\ref{teo:regularidad-eliptica-local-operadores-haces} then give
$v\in H^{r+2}(M)$, using a finite cover of $M$.
Induction yields all even orders and, by the
Sobolev--Morrey embeddings, $v\in C^\infty(M)$. The substitution
$v=(4\pi\tau)^{-m/4}e^{-f/2}$ gives
\eqref{eq:per-euler-mu-f}.

Diffeomorphism invariance follows from the change of variables formula.
Under $\mathbf{g}\mapsto c\mathbf{g}$, we have
$R\mapsto c^{-1}R$, $|\nabla f|^2\mapsto c^{-1}|\nabla f|^2$, and
$d\lambda\mapsto c^{m/2}d\lambda$; the factor
$(4\pi c\tau)^{-m/2}$ compensates for the change in volume. This proves
\eqref{eq:per-invariancias-mu}.
\end{proof}

If $M$ has connected components $M_1,\ldots,M_N$, set
$\displaystyle a_\alpha=\int_{M_\alpha}v^2\,d\lambda_{\mathbf{g}}$ for
$\alpha\in\{1,\ldots,N\}$. On the components with $a_\alpha>0$,
write $v=\sqrt{a_\alpha}\,w_\alpha$, with
$\|w_\alpha\|_{L^2(M_\alpha)}=1$. Then
\[
 \mathcal J_{\mathbf{g},\tau}(v)
 =\sum_{\substack{\alpha\in\{1,\ldots,N\}\\a_\alpha>0}}
       a_\alpha\mathcal J_{\mathbf{g}|_{M_\alpha},\tau}(w_\alpha)
  -\sum_{\alpha=1}^{N}a_\alpha\log a_\alpha,
 \qquad \sum_{\alpha=1}^{N}a_\alpha=1.
\]
The last term is nonnegative. Concentrating the mass on a component
that minimizes its quantity $\mu$, and approximating if positivity on
all components is required, gives
$\displaystyle \mu(\mathbf{g},\tau)=
\min_{\alpha\in\{1,\ldots,N\}}\mu(\mathbf{g}|_{M_\alpha},\tau)$.
Thus the following arguments apply component by component.

\begin{corollary}[Monotonicity of $\mu$]
\label{cor:per-monotonia-mu}
Let $M$ be a smooth closed connected manifold. If $\mathbf{g}(t)$ solves the
Ricci flow, $\tau'=-1$, and $\tau(t)>0$, then
$t\mapsto\mu(\mathbf{g}(t),\tau(t))$ is nondecreasing. If it remains
constant on an interval, the flow is, modulo diffeomorphisms and
rescaling, a gradient shrinking soliton.
\end{corollary}

\begin{proof}
Let $t_1<t_2$ and choose a positive minimizer for
$(\mathbf{g}(t_2),\tau(t_2))$. Its terminal density determines a positive
solution of the conjugate equation backward to $t_1$; by conservation
of mass, it can be written as
$(4\pi\tau)^{-m/2}e^{-f}$.
Theorem~\ref{teo:per-monotonia-W} gives
\[
 \mu(\mathbf{g}(t_2),\tau(t_2))=\mathcal W(t_2)
 \geq\mathcal W(t_1)\geq\mu(\mathbf{g}(t_1),\tau(t_1)).
\]
If the endpoints are equal, both inequalities are equalities, and the
square in \eqref{eq:per-monotonia-W} vanishes throughout the interval. This
is the gradient shrinking soliton equation.
\end{proof}

\section{Reduced length, distance, and volume}
\label{sec:per-geometria-reducida}

The entropy $\mathcal W$ yields noncollapsing for smooth flows. To
study the spacetime geometry more precisely, Perelman
introduced a second monotone quantity based on a variational problem for
curves. Its integrand combines kinetic energy with scalar curvature and
uses backward time as its parameter.

Throughout this section, $M$ is closed and connected. Fix $p\in M$, a time
$t_0$, and $\tau_{\max}>0$ such that
$[t_0-\tau_{\max},t_0]$ is contained in the smooth time interval of the flow. For
$0\leq s\leq\tau_{\max}$, write
\[
 \mathbf{g}_s:=\mathbf{g}(t_0-s),\qquad
 R_s:=R_{\mathbf{g}_s},\qquad
 \operatorname{Ric}_s:=\operatorname{Ric}_{\mathbf{g}_s}.
\]
We denote the Levi--Civita connection and
Laplacian of $\mathbf{g}_s$ by $\nabla^{(s)}$ and $\Delta_s$. Since $s$ is backward time,
$\partial_s\mathbf{g}_s=2\operatorname{Ric}_s$.

To differentiate a field along a curve, we use at each instant
the connection corresponding to that time. If $\mathbf{X}(s)=\dot\gamma(s)$ and
$\mathbf{Z}(s)\in T_{\gamma(s)}M$, we write
$\nabla_{\mathbf{X}}\mathbf{Z}=D_s\mathbf{Z}$ for the total covariant
derivative. In coordinates,
\begin{equation}
 (D_s\mathbf{Z})^k
 =\frac{dZ^k}{ds}
  +\sum_{i,j=1}^{m}\Gamma^k_{ij}(s,\gamma(s))X^iZ^j,
 \qquad k\in\{1,\ldots,m\}.
 \label{eq:per-derivada-total-curva}
\end{equation}
In particular, $D_s^2\mathbf{Z}$ includes the time derivative of the
Christoffel symbols appearing in $D_s\mathbf{Z}$. For two fields
$\mathbf{Z},\mathbf{W}$ along the curve, metric compatibility
on each slice gives
\begin{equation}
 \frac d{ds}\langle\mathbf{Z},\mathbf{W}\rangle_{\mathbf{g}_s}
 =\langle D_s\mathbf{Z},\mathbf{W}\rangle_{\mathbf{g}_s}
  +\langle\mathbf{Z},D_s\mathbf{W}\rangle_{\mathbf{g}_s}
  +2\operatorname{Ric}_s(\mathbf{Z},\mathbf{W}).
 \label{eq:per-producto-metrica-variable}
\end{equation}
Tensor derivatives involving only a spatial direction, such as
$\nabla_{\mathbf{Y}}\operatorname{Ric}_s$, are taken with the connection
$\nabla^{(s)}$ at fixed time $s$.

For an absolutely continuous curve
$\gamma\colon[0,\bar\tau]\to M$ with $\gamma(0)=p$ and
$\gamma(\bar\tau)=q$, define
\begin{equation}
 \mathcal L(\gamma)
 :=\int_0^{\bar\tau}\sqrt s\,
 \bigl(R_s(\gamma(s))+|\dot\gamma(s)|_{\mathbf{g}_s}^2\bigr)\,ds.
 \label{eq:per-longitud-L}
\end{equation}
The $\mathcal L$-distance and the reduced distance based at $(p,t_0)$
are
\begin{equation}
 L(q,\bar\tau):=\inf_{\gamma\in\mathcal A_{p,q}(\bar\tau)}\mathcal L(\gamma),
 \qquad l(q,\bar\tau):=\frac{L(q,\bar\tau)}{2\sqrt{\bar\tau}}.
 \label{eq:per-distancia-reducida}
\end{equation}
Here $\mathcal A_{p,q}(\bar\tau)$ is the set of absolutely continuous curves
from $[0,\bar\tau]$ into $M$ with these endpoints and finite weighted
energy, that is,
$\displaystyle \int_0^{\bar\tau}\sqrt s\,
|\dot\gamma(s)|_{\mathbf{g}_s}^2\,ds<\infty$.

The apparent singularity of the integrand at $s=0$ disappears after
reparametrization by the square root of time. This observation allows us to
apply the variational method and the theory of ordinary differential
equations directly.

\begin{lemma}[Regular reparametrization and existence of minimizers]
\label{lem:per-reparametrizacion-minimizantes-L}
Let $M$ be closed. If $\alpha(r)=\gamma(r^2)$ and
$\mathbf{A}(r)=\alpha'(r)$, then
\begin{equation}
 \mathcal L(\gamma)=\int_0^{\sqrt{\bar\tau}}
 \left(\frac12|\mathbf{A}(r)|_{\mathbf{g}_{r^2}}^2
 +2r^2R_{r^2}(\alpha(r))\right)dr.
 \label{eq:per-L-reparametrizada}
\end{equation}
For every $q\in M$ and every $0<\bar\tau\leq\tau_{\max}$, there exists an
$\mathcal L$-minimizing curve from $(p,0)$ to $(q,\bar\tau)$. Every minimizer is
smooth in the variable $r$ and satisfies
\begin{equation}
 \nabla^{(r^2)}_{\mathbf{A}}\mathbf{A}
 -2r^2\nabla^{\mathbf{g}_{r^2}}R_{r^2}
 +4r\operatorname{Ric}_{r^2}(\mathbf{A},\,\cdot\,)^{\sharp_{\mathbf{g}_{r^2}}}=0.
 \label{eq:per-geodesica-L-regular}
\end{equation}
In the original variable, with $\mathbf{X}=\dot\gamma$, this equation becomes
\begin{equation}
 \nabla^{(s)}_{\mathbf{X}}\mathbf{X}
 -\frac12\nabla^{\mathbf{g}_s}R_s+\frac1{2s}\mathbf{X}
 +2\operatorname{Ric}_s(\mathbf{X},\,\cdot\,)^{\sharp_{\mathbf{g}_s}}=0.
 \label{eq:per-geodesica-L}
\end{equation}
\end{lemma}

\begin{proof}
The change of variables $s=r^2$ gives $ds=2r\,dr$ and
$\mathbf{X}(r^2)=\mathbf{A}(r)/(2r)$. Direct substitution into
\eqref{eq:per-longitud-L} yields
\eqref{eq:per-L-reparametrizada}.

Fix a reference metric $\mathbf{h}$ on $M$. On the compact
interval under consideration, the metrics $\mathbf{g}_{r^2}$ are uniformly
equivalent to $\mathbf{h}$, and $R_{r^2}$ is uniformly bounded. Thus
there exist $c,C>0$ such that
\begin{equation}
 \mathcal L(\gamma)
 \geq c\int_0^{\sqrt{\bar\tau}}|\alpha'(r)|_{\mathbf{h}}^2\,dr-C.
 \label{eq:per-coercividad-L-reparametrizada}
\end{equation}
Choose a smooth embedding $\iota\colon M\hookrightarrow\mathbb R^N$.
The preceding estimate and the uniform equivalence of the metric induced
by $\iota$ and $\mathbf{h}$ bound $\iota\circ\alpha_j$ in
$H^1([0,\sqrt{\bar\tau}],\mathbb R^N)$ for every minimizing sequence.
After passing to a subsequence, these curves converge weakly in
$H^1([0,\sqrt{\bar\tau}],\mathbb R^N)$ and uniformly to a curve $a$. Since $\iota(M)$ is closed and each
$\iota\circ\alpha_j$ takes values in $\iota(M)$, so does $a$;
write $a=\iota\circ\alpha$, with values in $\iota(M)$. The endpoints are preserved
by uniform convergence.

Choose a finite partition $0=r_0<r_1<\cdots<r_J=\sqrt{\bar\tau}$
such that the $\alpha$-image of each interval
$I_\nu=[r_{\nu-1},r_\nu]$ lies in a compact subset of a chart
$(U_\nu,\varphi_\nu)$. Uniform convergence gives the same property for
$\alpha_j$ when $j$ is sufficiently large. The curves
$x_j=\varphi_\nu\circ\alpha_j$ converge uniformly to
$x=\varphi_\nu\circ\alpha$ and weakly in $H^1(I_\nu,\mathbb R^m)$.
Indeed, the chain rule and the bounds on the chart derivatives bound
their $H^1$ norms, while the uniform limit identifies every weak limit
with $x$.

Let $G_j(r)$ and $G(r)$ be the matrices of $\mathbf g_{r^2}$ evaluated
at $x_j(r)$ and $x(r)$, respectively. Then $G_j\to G$ uniformly on
$I_\nu$, and
\[
 \left|\int_{I_\nu}
       \langle(G_j-G)x_j',x_j'\rangle\,dr\right|
 \leq\|G_j-G\|_{L^\infty(I_\nu)}
       \|x_j'\|_{L^2(I_\nu,\mathbb R^m)}^2\longrightarrow0.
\]
Multiplication by the continuous matrix $G^{1/2}(r)$ is bounded on
$L^2(I_\nu,\mathbb R^m)$. Hence
$G^{1/2}x_j'\rightharpoonup G^{1/2}x'$ in this space, and weak lower
semicontinuity of the norm yields
\[
 \int_{I_\nu}\langle Gx',x'\rangle\,dr
 \leq\liminf_{j\to\infty}
       \int_{I_\nu}\langle G_jx_j',x_j'\rangle\,dr.
\]
Summing over the finite partition gives
\[
 \int_0^{\sqrt{\bar\tau}}|\alpha'|_{\mathbf g_{r^2}}^2\,dr
 \leq\liminf_{j\to\infty}
       \int_0^{\sqrt{\bar\tau}}|\alpha_j'|_{\mathbf g_{r^2}}^2\,dr.
\]
The curvature term converges by uniform convergence and the dominated
convergence theorem. Thus the action of $\alpha$ does not exceed the
limit of the actions of the minimizing sequence.
The curve $\gamma(s)=\alpha(\sqrt s)$ is admissible: the substitution
$s=r^2$ gives
\[
 \int_0^{\bar\tau}|\dot\gamma(s)|_{\mathbf h}\,ds
 =\int_0^{\sqrt{\bar\tau}}|\alpha'(r)|_{\mathbf h}\,dr<\infty,
 \qquad
 \int_0^{\bar\tau}\sqrt s\,|\dot\gamma(s)|_{\mathbf g_s}^2\,ds
 =\frac12\int_0^{\sqrt{\bar\tau}}|\alpha'(r)|_{\mathbf g_{r^2}}^2\,dr.
\]
The first equality, also applied on subintervals, proves absolute
continuity up to $s=0$; the second proves finiteness of the weighted
energy. Therefore $\gamma$ attains the infimum of $\mathcal L$.

The first variation of \eqref{eq:per-L-reparametrizada} gives the
Euler--Lagrange equation in the weak sense. To justify its regularity,
consider a subinterval $I\subseteq[0,\sqrt{\bar\tau}]$ whose image lies in a chart and define
the momenta
\[
 P_i(r)=\sum_{j=1}^{m}g_{ij}(r^2,\alpha(r))\alpha'^j(r),
 \qquad i\in\{1,\ldots,m\}.
\]
The weak equation is
\[
 P_i'(r)=\frac12\sum_{j,k=1}^{m}
  \partial_i g_{jk}(r^2,\alpha(r))\alpha'^j(r)\alpha'^k(r)
  +2r^2\partial_iR_{r^2}(\alpha(r)).
\]
The right-hand side belongs to $L^1(I)$ because $\alpha'\in L^2(I)$.
Thus each $P_i$ has an absolutely continuous, and hence continuous, representative.
The identity
$\displaystyle \alpha'^j=\displaystyle\sum_{i=1}^{m}g^{ji}(r^2,\alpha)P_i$
then shows that $\alpha$ is $C^1$. The right-hand side of the equation
for $P_i'$ is now continuous, so $P_i$ is
$C^1$ and $\alpha$ is $C^2$. The resulting ordinary differential system has smooth coefficients, even
at $r=0$, and ordinary differential equation regularity gives smoothness up to the endpoints.
Writing it with the connection of $\mathbf{g}_{r^2}$ gives
\eqref{eq:per-geodesica-L-regular}.
Finally,
$\mathbf{A}=2\sqrt s\,\mathbf{X}$ and
\[
 \nabla^{(r^2)}_{\mathbf{A}}\mathbf{A}
 =2\mathbf{X}+4r^2\nabla^{(s)}_{\mathbf{X}}\mathbf{X}.
\]
Substituting this identity into
\eqref{eq:per-geodesica-L-regular} proves
\eqref{eq:per-geodesica-L}.
\end{proof}

\begin{lemma}[The $\mathcal L$-exponential map]
\label{lem:per-exponencial-L}
For each $v\in T_pM$, there exists a unique solution of
\eqref{eq:per-geodesica-L-regular} with
\[
 \alpha_v(0)=p,
 \qquad \alpha_v'(0)=2v.
\]
It is defined on all of $[0,\sqrt{\tau_{\max}}]$ and depends smoothly on
$(v,r)$. If $\gamma_v(s)=\alpha_v(\sqrt s)$, then
\[
 \lim_{s\to 0^{+}}\sqrt s\,\dot\gamma_v(s)=v.
\]
Define
\begin{equation}
 \mathcal L\exp_p^\tau(v):=\gamma_v(\tau).
 \label{eq:per-exponencial-L}
\end{equation}
\end{lemma}

\begin{proof}
Existence, uniqueness, and smooth dependence near $r=0$ follow from the
usual theorem for ordinary differential systems applied to
\eqref{eq:per-geodesica-L-regular}. To exclude blowup of the
velocity on a finite interval, take the inner product of that equation with
$\mathbf{A}$. The uniform bounds on $\nabla R$ and $\operatorname{Ric}$
give
\[
 \frac d{dr}|\mathbf{A}|_{\mathbf{g}_{r^2}}^2
 \leq C\bigl(1+|\mathbf{A}|_{\mathbf{g}_{r^2}}^2\bigr).
\]
Grönwall's inequality keeps the velocity bounded. Since $M$ is
compact, the solution extends to
$\sqrt{\tau_{\max}}$. The identity
$\mathbf{A}(r)=2r\mathbf{X}(r^2)$ gives the last assertion.
\end{proof}

\begin{definition}[$\mathcal L$-injectivity domain and cut locus]
\label{def:per-lugar-corte-L}
For $\tau>0$, let $\widetilde U_p(\tau)\subseteq T_pM$ be the set of
vectors $v$ satisfying the following two conditions:
\begin{enumerate}
 \item $d(\mathcal L\exp_p^\tau)_v$ is an isomorphism;
 \item there exists a neighborhood $\mathcal V$ of $v$ in $T_pM$ such that, for every
 $w\in\mathcal V$, the curve $\gamma_w|_{[0,\tau]}$ is the unique
 $\mathcal L$-minimizer between its endpoints.
\end{enumerate}
Set
\[
 U_p(\tau):=\mathcal L\exp_p^\tau(\widetilde U_p(\tau)),
 \qquad \operatorname{Cut}_{\mathcal L}(p,\tau):=M\setminus U_p(\tau).
\]
The second set is called the \emph{$\mathcal L$-cut locus} based
at $(p,t_0)$ at backward time $\tau$.
\end{definition}

By definition, $\widetilde U_p(\tau)$ is open. The restriction of
$\mathcal L\exp_p^\tau$ to this set is injective: two preimages of the
same point would determine two minimizers, but each is unique. The
inverse function theorem therefore shows that this restriction is a
diffeomorphism onto the open set $U_p(\tau)$. In particular, $L(\cdot,\tau)$ and
$l(\cdot,\tau)$ are smooth on $U_p(\tau)$.

\begin{lemma}[First variation]
\label{lem:per-primera-variacion-L}
Let $\gamma_u$ be a variation of an $\mathcal L$--geodesic
$\gamma=\gamma_0$, smooth in $(u,r)$ with $r=\sqrt s$, that fixes the initial endpoint
$\gamma_u(0)=p$, and let
$\mathbf{Y}=\left.\partial_u\gamma_u\right|_{u=0}$. Then
$\mathbf{Y}(0)=0$ and
\begin{align}
 \delta\mathcal L(\gamma)
 =2\sqrt s\,\langle\mathbf{Y},\mathbf{X}\rangle_{\mathbf{g}_s}
 \Big|_{0}^{\bar\tau}+\int_0^{\bar\tau}\sqrt s\left\langle\mathbf{Y},
 \nabla R_s-\frac1s\mathbf{X}-2\nabla_{\mathbf{X}}\mathbf{X}-4\operatorname{Ric}_s(\mathbf{X},\,\cdot\,)^{\sharp}
 \right\rangle_{\mathbf{g}_s}ds.
 \label{eq:per-primera-variacion-L-completa}
\end{align}
In particular, if $q\in U_p(\tau)$ and $\gamma$ is the regular minimizer
ending at $q$, then
\begin{equation}
 \nabla^{\mathbf{g}_\tau}L(q,\tau)=2\sqrt\tau\,\mathbf{X}(\tau),
 \qquad
 \partial_\tau L(q,\tau)=\sqrt\tau
 \bigl(R_\tau(q)-|\mathbf{X}(\tau)|_{\mathbf{g}_\tau}^2\bigr).
 \label{eq:per-primera-variacion-L}
\end{equation}
Consequently,
\begin{equation}
 \nabla^{\mathbf{g}_\tau}l=\mathbf{X}(\tau),
 \qquad
 \partial_\tau l=-\frac l{2\tau}
 +\frac12\bigl(R_\tau-|\nabla l|_{\mathbf{g}_\tau}^2\bigr)
 \label{eq:per-primera-variacion-l}
\end{equation}
on the regular set.
\end{lemma}

\begin{proof}
Let $\gamma_u$ be a variation and let $\mathbf{Y}=\partial_u\gamma_u|_{u=0}$.
Since $\partial_s\mathbf{g}_s=2\operatorname{Ric}_s$ and the coordinate fields
of the variation commute,
\[
 \partial_u|\mathbf{X}|_{\mathbf{g}_s}^2
 =2\langle\nabla_{\mathbf{X}}\mathbf{Y},\mathbf{X}\rangle_{\mathbf{g}_s}.
\]
When integrating by parts, we must also differentiate the
$s$-dependent metric:
\begin{align*}
 \frac d{ds}\bigl(2\sqrt s\langle\mathbf{Y},\mathbf{X}\rangle\bigr)
 =\frac1{\sqrt s}\langle\mathbf{Y},\mathbf{X}\rangle
 +2\sqrt s\langle\nabla_{\mathbf{X}}\mathbf{Y},\mathbf{X}\rangle+2\sqrt s\langle\mathbf{Y},\nabla_{\mathbf{X}}\mathbf{X}\rangle
 +4\sqrt s\operatorname{Ric}_s(\mathbf{Y},\mathbf{X}).
\end{align*}
Substituting this identity into the derivative of
\eqref{eq:per-longitud-L} gives
\eqref{eq:per-primera-variacion-L-completa}. The endpoint term at
$s=0$ vanishes because $\mathbf{Y}(0)=0$ and
$\sqrt s\,\mathbf{X}(s)$ has a limit.

On the regular set, the inverse function theorem applied to
$\mathcal L\exp_p^\tau$ produces a smooth family of minimizers as the
endpoint varies. The boundary part of
\eqref{eq:per-primera-variacion-L-completa} gives the first identity in
\eqref{eq:per-primera-variacion-L}. Along the minimizer itself,
\[
 \frac d{d\tau}L(\gamma(\tau),\tau)
 =\sqrt\tau(R_\tau+|\mathbf{X}(\tau)|^2).
\]
By the chain rule, the left-hand side is also
$\partial_\tau L+\langle\nabla L,\mathbf{X}(\tau)\rangle$. The gradient identity gives
the second formula. Dividing by $2\sqrt\tau$ proves
\eqref{eq:per-primera-variacion-l}.
\end{proof}

\begin{lemma}[$\mathcal L$-Jacobi fields and index]
\label{lem:per-jacobi-indice-L}
Let $\mathbf{C}_s=\partial_s\nabla^{(s)}$ be the variation tensor of the
connection. Along an $\mathcal L$--geodesic, the linearization of
\eqref{eq:per-geodesica-L}, with the total derivative in
\eqref{eq:per-derivada-total-curva}, is
\begin{align}
 \mathcal J_{\mathcal L}\mathbf{Y}:={}&
 D_s^2\mathbf{Y}+\operatorname{Rm}_{\mathbf{g}_s}(\mathbf{Y},\mathbf{X})\mathbf{X}
 -\mathbf{C}_s(\mathbf{Y},\mathbf{X})
 -\frac12\nabla_{\mathbf{Y}}(\nabla R_s)
 +\frac1{2s}D_s\mathbf{Y}\nonumber\\
 &+2(\nabla_{\mathbf{Y}}\operatorname{Ric}_s)
       (\mathbf{X},\,\cdot\,)^{\sharp}
  +2\operatorname{Ric}_s(D_s\mathbf{Y},\,\cdot\,)^{\sharp}.
 \label{eq:per-operador-jacobi-L}
\end{align}
The fields satisfying $\mathcal J_{\mathcal L}\mathbf{Y}=0$ that are smooth
in $r=\sqrt s$ and vanish at $s=0$ are precisely the derivatives of
families of $\mathcal L$--geodesics with fixed initial point $p$.

For fields $\mathbf{Y},\mathbf{Z}$ smooth in $r=\sqrt s$ and vanishing at
$s=0$, define
\begin{align}
 I_{\mathcal L}(\mathbf{Y},\mathbf{Z}):={}&
 \int_0^\tau\sqrt s\,\Bigl\{
 \operatorname{Hess}_{\mathbf{g}_s}R_s(\mathbf{Y},\mathbf{Z})
 +2\langle\operatorname{Rm}_{\mathbf{g}_s}(\mathbf{Y},\mathbf{X})
               \mathbf{Z},\mathbf{X}\rangle_{\mathbf{g}_s}\nonumber\\
 &+2\langle D_s\mathbf{Y},D_s\mathbf{Z}\rangle_{\mathbf{g}_s}
 +2(\nabla_{\mathbf{X}}\operatorname{Ric}_s)(\mathbf{Y},\mathbf{Z})
 -2(\nabla_{\mathbf{Y}}\operatorname{Ric}_s)(\mathbf{X},\mathbf{Z})\nonumber\\
 &-2(\nabla_{\mathbf{Z}}\operatorname{Ric}_s)(\mathbf{X},\mathbf{Y})
 \Bigr\}\,ds.
 \label{eq:per-forma-indice-L}
\end{align}
If $\gamma_u$ is a variation of a minimizer, smooth in $(u,\sqrt s)$,
with fixed initial point and variation field $\mathbf{Y}$, then
\begin{equation}
 \left.\frac{d^2}{du^2}\right|_{u=0}\mathcal L(\gamma_u)
 =2\sqrt\tau\,
 \langle\nabla_{\mathbf{Y}}\mathbf{Y},\mathbf{X}\rangle_{\mathbf{g}_\tau}
 \Big|_{s=\tau}+I_{\mathcal L}(\mathbf{Y},\mathbf{Y}).
 \label{eq:per-segunda-variacion-L}
\end{equation}
At a regular point,
\begin{equation}
 \operatorname{Hess}_{\mathbf{g}_\tau}L
 (\mathbf{Y}(\tau),\mathbf{Y}(\tau))
 \leq I_{\mathcal L}(\mathbf{Y},\mathbf{Y}).
 \label{eq:per-cota-hessiana-indice-L}
\end{equation}
Moreover, $d(\mathcal L\exp_p^\tau)_v$ is singular if and only if there exists a
nonzero Jacobi field vanishing at both endpoints.
\end{lemma}

\begin{proof}
The connection variation formula from the previous chapter, applied to
$\partial_s\mathbf{g}_s=2\operatorname{Ric}_s$, gives
\begin{align}
 \langle\mathbf{C}_s(\mathbf{U},\mathbf{V}),\mathbf{W}\rangle_{\mathbf{g}_s}
 ={}&(\nabla_{\mathbf{U}}\operatorname{Ric}_s)(\mathbf{V},\mathbf{W})
 +(\nabla_{\mathbf{V}}\operatorname{Ric}_s)(\mathbf{U},\mathbf{W})
 -(\nabla_{\mathbf{W}}\operatorname{Ric}_s)(\mathbf{U},\mathbf{V}).
 \label{eq:per-variacion-conexion-retrograda}
\end{align}
Here $\mathbf{U},\mathbf{V},\mathbf{W}$ are tangent vectors at the same
point, and all derivatives on the right-hand side are spatial.

On a variation $\gamma(u,s)$, write $D_u$ for the covariant
derivative in the $\mathbf{Y}=\partial_u\gamma$ direction at fixed $s$.
The absence of torsion gives $D_u\mathbf{X}=D_s\mathbf{Y}$. For a field
$\mathbf{Z}$ on the variation, the commutator is
\begin{equation}
 D_uD_s\mathbf{Z}-D_sD_u\mathbf{Z}
 =\operatorname{Rm}_{\mathbf{g}_s}(\mathbf{Y},\mathbf{X})\mathbf{Z}
  -\mathbf{C}_s(\mathbf{Y},\mathbf{Z}).
 \label{eq:per-conmutador-variacion-tiempo}
\end{equation}
To check the time-dependent term, write each derivative as in
\eqref{eq:per-derivada-total-curva}. The terms containing spatial
derivatives of $\Gamma^k_{ij}$ and products of two Christoffel symbols
combine into the curvature. Differentiating $D_u\mathbf{Z}$ with respect to $s$
also produces
$\displaystyle \sum_{i,j=1}^{m}(\partial_s\Gamma^k_{ij})Y^iZ^j$,
for each $k\in\{1,\ldots,m\}$, with a negative sign in the commutator.
This is exactly the last summand in
\eqref{eq:per-conmutador-variacion-tiempo}.

In particular,
\[
 D_u(D_s\mathbf{X})
 =D_s^2\mathbf{Y}+\operatorname{Rm}(\mathbf{Y},\mathbf{X})\mathbf{X}
  -\mathbf{C}_s(\mathbf{Y},\mathbf{X}).
\]
The derivative of $\operatorname{Ric}_s(\mathbf{X},\cdot)^{\sharp}$ in the
$u$ direction is
$(\nabla_{\mathbf{Y}}\operatorname{Ric}_s)(\mathbf{X},\cdot)^{\sharp}
+\operatorname{Ric}_s(D_s\mathbf{Y},\cdot)^{\sharp}$.
Differentiating the remaining terms of \eqref{eq:per-geodesica-L} yields
\eqref{eq:per-operador-jacobi-L}.

The equation in $r=\sqrt s$ is a regular second-order linear ordinary
differential system. If $\mathbf{Y}(0)=0$, it is determined by
$\displaystyle \lim_{s\to0^{+}}\sqrt s\,D_s\mathbf{Y}$.
Indeed, for a family of initial data $v(u)$,
$D_r\mathbf{Y}(0)=2v'(0)$, and this limit is $v'(0)$.
Smooth dependence and uniqueness of the initial value problem establish the
correspondence with variations of $\mathcal L$--geodesics.

Now differentiate \eqref{eq:per-primera-variacion-L-completa} along the
central curve. Its Euler--Lagrange equation makes the integrand of the
first variation vanish. We obtain
\begin{align*}
 \delta^2\mathcal L
 ={}&2\sqrt\tau\bigl(
  \langle\nabla_{\mathbf{Y}}\mathbf{Y},\mathbf{X}\rangle
  +\langle D_s\mathbf{Y},\mathbf{Y}\rangle\bigr)_{s=\tau}
 -2\int_0^\tau\sqrt s\,
       \langle\mathcal J_{\mathcal L}\mathbf{Y},\mathbf{Y}\rangle\,ds.
\end{align*}
To integrate the last term, use
\eqref{eq:per-producto-metrica-variable} with
$\mathbf{Z}=D_s\mathbf{Y}$ and $\mathbf{W}=\mathbf{Y}$:
\begin{equation}
 \frac d{ds}\langle D_s\mathbf{Y},\mathbf{Y}\rangle_{\mathbf{g}_s}
 =\langle D_s^2\mathbf{Y},\mathbf{Y}\rangle_{\mathbf{g}_s}
  +|D_s\mathbf{Y}|_{\mathbf{g}_s}^2
  +2\operatorname{Ric}_s(D_s\mathbf{Y},\mathbf{Y}).
 \label{eq:per-derivada-metrica-indice-L}
\end{equation}
On the other hand, setting $\mathbf{U}=\mathbf{W}=\mathbf{Y}$ and
$\mathbf{V}=\mathbf{X}$ in \eqref{eq:per-variacion-conexion-retrograda},
the two terms involving $\nabla_{\mathbf{Y}}\operatorname{Ric}_s$ cancel:
\[
 \langle\mathbf{C}_s(\mathbf{Y},\mathbf{X}),\mathbf{Y}\rangle
 =(\nabla_{\mathbf{X}}\operatorname{Ric}_s)(\mathbf{Y},\mathbf{Y}).
\]
Consequently,
\begin{align*}
 -2\sqrt s\,\langle\mathcal J_{\mathcal L}\mathbf{Y},\mathbf{Y}\rangle
 ={}&-\frac d{ds}\bigl(2\sqrt s\langle D_s\mathbf{Y},\mathbf{Y}\rangle\bigr)\\
 &+\sqrt s\Bigl\{\operatorname{Hess}R_s(\mathbf{Y},\mathbf{Y})
 +2\langle\operatorname{Rm}(\mathbf{Y},\mathbf{X})\mathbf{Y},\mathbf{X}\rangle
 +2|D_s\mathbf{Y}|^2\\
 &\hspace{22mm}+2(\nabla_{\mathbf{X}}\operatorname{Ric}_s)(\mathbf{Y},\mathbf{Y})
 -4(\nabla_{\mathbf{Y}}\operatorname{Ric}_s)(\mathbf{X},\mathbf{Y})\Bigr\}.
\end{align*}
The total derivative term cancels the summand
$2\sqrt\tau\langle D_s\mathbf{Y},\mathbf{Y}\rangle$ at the final endpoint.
At the initial endpoint, smoothness in $r$ and $\mathbf{Y}(0)=0$ imply
$\mathbf{Y}=O(\sqrt s)$ and $D_s\mathbf{Y}=O(s^{-1/2})$; hence
$2\sqrt s\langle D_s\mathbf{Y},\mathbf{Y}\rangle\to0$.
These same bounds justify convergence of the improper integrals.
This proves \eqref{eq:per-segunda-variacion-L}, and polarization gives
\eqref{eq:per-forma-indice-L}. If both fields also vanish at
$s=\tau$, bilinear integration by parts gives
\[
 I_{\mathcal L}(\mathbf{Y},\mathbf{Z})
 =-2\int_0^\tau\sqrt s\,
       \langle\mathcal J_{\mathcal L}\mathbf{Y},\mathbf{Z}\rangle\,ds.
\]
The symmetry of this expression also follows directly from
\eqref{eq:per-forma-indice-L} and the curvature symmetries.

To estimate the Hessian, realize $\mathbf{Y}$ by a variation
with fixed initial point and final endpoint $q(u)$. We have
$L(q(u),\tau)\leq\mathcal L(\gamma_u)$, with equality at $u=0$.
The second derivative of the difference is nonnegative. When we use
$\nabla L=2\sqrt\tau\,\mathbf{X}(\tau)$, the term containing the
acceleration of $q(u)$ is the same on both sides and cancels. This
proves \eqref{eq:per-cota-hessiana-indice-L}.
Finally, $d(\mathcal L\exp_p^\tau)_v(w)=\mathbf{Y}_w(\tau)$, where
$\displaystyle \lim_{s\to0^{+
}}\sqrt s D_s\mathbf{Y}_w=w$. Uniqueness for the regular system
implies that $w\ne0$ corresponds to a nonzero field, which proves
the assertion about the differential of the exponential map.
\end{proof}

\begin{lemma}[Identities for $K$ and a Hessian estimate]
\label{lem:per-identidades-K-Hessiana-L}
Let $q\in U_p(\tau)$ and let $\gamma$ be the regular minimizer ending at
$q$. Along $\gamma$, define
\begin{align}
 \mathcal H(\mathbf{X})
 &:=-\partial_sR_s-\frac{R_s}{s}
 -2\langle\nabla R_s,\mathbf{X}\rangle
 +2\operatorname{Ric}_s(\mathbf{X},\mathbf{X}),
 \label{eq:per-Hamilton-integrando-L}\\
 K(\tau)&:=\int_0^\tau s^{3/2}\mathcal H(\mathbf{X}(s))\,ds.
 \label{eq:per-K-longitud-reducida}
\end{align}
Then
\begin{align}
 \partial_\tau l
 &=R_\tau-\frac l\tau+\frac{K(\tau)}{2\tau^{3/2}},
 \label{eq:per-tiempo-l-K}\\
 |\nabla^{\mathbf{g}_\tau}l|_{\mathbf{g}_\tau}^2
 &=-R_\tau+\frac l\tau-\frac{K(\tau)}{\tau^{3/2}},
 \label{eq:per-gradiente-l-K}\\
 \Delta_\tau l
 &\leq\frac m{2\tau}-R_\tau-
 \frac{K(\tau)}{2\tau^{3/2}}.
 \label{eq:per-laplaciano-l-K}
\end{align}
Consequently, on the regular set,
\begin{equation}
 \partial_\tau l-\Delta_\tau l+|\nabla l|_{\mathbf{g}_\tau}^2-R_\tau
 +\frac m{2\tau}\geq0.
 \label{eq:per-desigualdad-l-reducida-regular}
\end{equation}
\end{lemma}

\begin{proof}
Differentiate $R_s(\gamma(s))+|\mathbf{X}(s)|^2$. The time variation of the
metric contributes $2\operatorname{Ric}_s(\mathbf{X},\mathbf{X})$, and equation
\eqref{eq:per-geodesica-L} gives
\begin{equation}
 \frac d{ds}\bigl(R_s+|\mathbf{X}|^2\bigr)
 =-\mathcal H(\mathbf{X})-\frac1s(R_s+|\mathbf{X}|^2).
 \label{eq:per-derivada-R-X-L}
\end{equation}
Thus
\[
 \frac d{ds}\left[s^{3/2}(R_s+|\mathbf{X}|^2)\right]
 =\frac12\sqrt s(R_s+|\mathbf{X}|^2)-s^{3/2}\mathcal H(\mathbf{X}).
\]
Integration between $0$ and $\tau$, together with
$\displaystyle \lim_{s\to 0^{+}}s^{3/2}(R_s+|\mathbf{X}|^2)=0$, yields
\begin{equation}
 \tau^{3/2}(R_\tau+|\mathbf{X}(\tau)|^2)
 =\frac12L(q,\tau)-K(\tau).
 \label{eq:per-identidad-basica-K}
\end{equation}
The second formula in \eqref{eq:per-primera-variacion-L} and the identity
$\nabla l=\mathbf{X}(\tau)$ transform
\eqref{eq:per-identidad-basica-K} into
\eqref{eq:per-tiempo-l-K} and \eqref{eq:per-gradiente-l-K}.

It remains to estimate the trace of the Hessian. Let
$(\mathbf{e}_a)_{a=1}^m$ be an orthonormal basis of $T_qM$, and let
$\mathbf{Y}_a$ be the comparison field determined by
\begin{equation}
 \mathbf{Y}_a(\tau)=\mathbf{e}_a,
 \qquad
 D_s\mathbf{Y}_a
 =-\operatorname{Ric}_s(\mathbf{Y}_a,\,\cdot\,)^{\sharp}
 +\frac1{2s}\mathbf{Y}_a.
 \label{eq:per-campos-comparacion-L}
\end{equation}
Differentiating
$\langle\mathbf{Y}_a,\mathbf{Y}_b\rangle_{\mathbf{g}_s}$ and using
$\partial_s\mathbf{g}_s=2\operatorname{Ric}_s$ gives
\[
 \langle\mathbf{Y}_a(s),\mathbf{Y}_b(s)\rangle_{\mathbf{g}_s}
 =\frac s\tau\delta_{ab}.
\]
Thus the fields vanish at $s=0$, and
$\mathbf{W}_a(s)=\sqrt{\tau/s}\,\mathbf{Y}_a(s)$ form an orthonormal frame
satisfying
$D_s\mathbf{W}_a
=-\operatorname{Ric}_s(\mathbf{W}_a,\cdot)^\sharp$ for
$a\in\{1,\ldots,m\}$. In the variable $r=\sqrt s$, this equation is
$D_r\mathbf{W}_a=-2r\operatorname{Ric}_{r^2}(\mathbf{W}_a,\cdot)^\sharp$;
its coefficients are smooth up to $r=0$. Hence the fields
$\mathbf{W}_a$ extend smoothly, and
$\mathbf{Y}_a=(r/\sqrt\tau)\mathbf{W}_a$ are admissible variations.
Inequality \eqref{eq:per-cota-hessiana-indice-L} implies
\[
 \Delta_\tau L(q,\tau)
 \leq\sum_{a=1}^mI_{\mathcal L}(\mathbf{Y}_a,\mathbf{Y}_a).
\]

Substitute $\mathbf{Y}_a=\sqrt{s/\tau}\,\mathbf{W}_a$ into
\eqref{eq:per-forma-indice-L}. The required contractions are
\begin{align*}
 \sum_{a=1}^{m}\operatorname{Hess}R_s(\mathbf{W}_a,\mathbf{W}_a)&=\Delta_sR_s,\\
 \sum_{a=1}^{m}\langle\operatorname{Rm}(\mathbf{W}_a,\mathbf{X})
 \mathbf{W}_a,\mathbf{X}\rangle&=-\operatorname{Ric}_s(\mathbf{X},\mathbf{X}),\\
 \sum_{a=1}^{m}\left|\frac1{2s}\mathbf{W}_a-
 \operatorname{Ric}_s(\mathbf{W}_a,\,\cdot\,)^{\sharp}\right|^2
 &=\frac m{4s^2}-\frac{R_s}{s}+|\operatorname{Ric}_s|^2,\\
 \sum_{a=1}^{m}(\nabla_{\mathbf{W}_a}\operatorname{Ric}_s)
 (\mathbf{X},\mathbf{W}_a)&=\frac12\langle\nabla R_s,\mathbf{X}\rangle.
\end{align*}
The last equality is the contracted Bianchi identity. Together with
$\partial_sR_s=-\Delta_sR_s-2|\operatorname{Ric}_s|^2$, these identities
give, term by term,
\begin{align}
 \sum_{a=1}^mI_{\mathcal L}(\mathbf{Y}_a,\mathbf{Y}_a)
 =\frac1\tau\int_0^\tau\left[
 \frac m{2\sqrt s}-2\sqrt sR_s-s^{3/2}\partial_sR_s
 -2s^{3/2}\operatorname{Ric}_s(\mathbf{X},\mathbf{X})\right]ds.
 \label{eq:per-traza-indice-L-explicita}
\end{align}
Here all quantities are evaluated at $\gamma(s)$. Now
\[
 \frac d{ds}\bigl(2s^{3/2}R_s(\gamma(s))\bigr)
 =3\sqrt sR_s+2s^{3/2}(\partial_sR_s+\langle\nabla R_s,\mathbf{X}\rangle).
\]
The definition of $\mathcal H$ allows us to rewrite the part of
the preceding integral that does not involve the dimension as
\[
 -2\tau^{3/2}R_\tau(q)-K(\tau).
\]
Since $\displaystyle \int_0^\tau m/(2\sqrt s)\,ds=m\sqrt\tau$, we conclude that
\[
 \Delta_\tau L(q,\tau)
 \leq\frac m{\sqrt\tau}-2\sqrt\tau R_\tau(q)-\frac{K(\tau)}\tau.
\]
Dividing by $2\sqrt\tau$ proves
\eqref{eq:per-laplaciano-l-K}. Combining the three formulas gives
\eqref{eq:per-desigualdad-l-reducida-regular}.
\end{proof}

The preceding inequality was obtained where the minimizer is unique and the
exponential map has no conjugate points. To integrate it over the entire
manifold, we need to control the reduced function on its cut locus.

\begin{definition}[Local semiconcavity]
\label{def:per-semiconcavidad-local}
\index{local semiconcavity}
Let $\Omega\subseteq\mathbb R^m$ be open. A continuous function
$u\colon\Omega\longrightarrow\mathbb R$ is \emph{locally semiconcave}
if each point has a convex neighborhood $V\subseteq\Omega$ on which,
for some constant $C\geq0$, the function
\[
 x\longmapsto u(x)-\frac C2|x|^2
\]
is concave. On a smooth manifold, this condition is imposed on the
coordinate representations of $u$. The constants may depend on the
chosen chart and neighborhood; no global constant is required.
\end{definition}

Concavity of a function $v$ on $V$ means that
\[
 v((1-t)x+ty)\geq(1-t)v(x)+tv(y),
 \qquad x,y\in V,\quad 0\leq t\leq1.
\]
For a continuous function, it is enough to check this inequality at
midpoints: repeated application of that inequality gives all
dyadic coefficients, and continuity allows passage to any
$t\in[0,1]$. We will use this observation both to recognize
semiconcavity and to check that the definition is independent of
coordinates.

\begin{lemma}[Upper supports and changes of coordinates]
\label{lem:per-apoyos-semiconcavidad}
Let $V\subseteq\mathbb R^m$ be open and convex, and let $u\in C(V)$.
Suppose that for each $x\in V$ there exists $\phi_x\in C^2(V)$ such that
$u\leq\phi_x$ on $V$, $u(x)=\phi_x(x)$, and
\[
 D^2\phi_x(y)[e,e]\leq C|e|^2,
 \qquad y\in V,\quad e\in\mathbb R^m,
\]
with the same constant $C\geq0$. Then $u-C|x|^2/2$ is concave on
$V$. Moreover, every locally semiconcave function is locally Lipschitz,
and local semiconcavity is preserved under composition with a
$C^2$ change of coordinates.
\end{lemma}

\begin{proof}
Let $x-h,x+h\in V$. Since $V$ is convex, the segment joining these points
lies in $V$. Taylor's formula with integral remainder, applied to
$t\mapsto\phi_x(x+th)$, gives
\[
 \phi_x(x+h)+\phi_x(x-h)-2\phi_x(x)
 =\int_{-1}^{1}(1-|t|)D^2\phi_x(x+th)[h,h]dt
 \leq C|h|^2.
\]
The properties of the upper support imply
\[
 u(x+h)+u(x-h)-2u(x)\leq C|h|^2.
\]
After subtracting $C|x|^2/2$, this inequality becomes midpoint
concavity. The preceding observation proves the first assertion.

Let us prove local Lipschitz regularity. It suffices to prove it for a
continuous concave function $v$, since the quadratic term is smooth. Choose
$\overline{B(x_0,2r)}\subset V$ and set
$M_r=\displaystyle\max_{\overline{B(x_0,2r)}}|v|$. If $x,y\in B(x_0,r)$ and $x\ne y$,
define $d=|y-x|$ and $e=(y-x)/d$. The points $x-re$ and $y+re$ belong
to $B(x_0,2r)$. Concavity along the line containing these four
points gives
\[
 \frac{v(x)-v(x-re)}r
 \geq\frac{v(y)-v(x)}d
 \geq\frac{v(y+re)-v(y)}r.
\]
Both outer quotients have absolute value at most $2M_r/r$.
Therefore, $|v(y)-v(x)|\leq(2M_r/r)|y-x|$.

Now let $\kappa$ be a $C^2$ change of coordinates. Shrink its
domains so that $u$ is $L$--Lipschitz on a convex open set
containing the image under consideration, $u-C|x|^2/2$ is concave there, and
$\|D\kappa\|\leq K_0$ and $\|D^2\kappa\|\leq K_1$. This is possible by
working in neighborhoods with compact closure contained in the domains.
For $x-h,x+h$ in a convex ball in the smaller domain, set
\[
 z_\pm=\kappa(x\pm h),\qquad z_0=\kappa(x),\qquad
 z_*=(z_++z_-)/2.
\]
Semiconcavity at the points $z_-,z_+$ and the Lipschitz bound give
\[
 u(z_+)+u(z_-)-2u(z_0)
 \leq\frac C4|z_+-z_-|^2+2L|z_*-z_0|.
\]
The mean value theorem and the same Taylor formula imply
\[
 |z_+-z_-|\leq2K_0|h|,\qquad
 |z_*-z_0|\leq\frac{K_1}{2}|h|^2.
\]
Consequently, the second increment of $u\circ\kappa$ is bounded
above by $(CK_0^2+LK_1)|h|^2$. The midpoint criterion
proves its semiconcavity. Applying the same argument to $\kappa^{-1}$ establishes
the chart independence asserted in the definition.
\end{proof}

\begin{lemma}[Regularity and cut locus of the reduced distance]
\label{lem:per-regularidad-lugar-corte-L}
For each compact interval
$[\tau_0,\tau_1]\subset(0,\tau_{\max}]$, the functions $L$ and $l$ are
locally Lipschitz on $M\times[\tau_0,\tau_1]$. For each $0<\tau\leq\tau_{\max}$,
$U_p(\tau)$ is open, $L$ and $l$ are smooth there, and
$\operatorname{Cut}_{\mathcal L}(p,\tau)$ has measure zero. Moreover, the
restriction $l(\cdot,\tau)$ is locally semiconcave.
\end{lemma}

\begin{proof}
Let $0<\tau_0\leq\tau_1\leq\tau_{\max}$ and put
$a_0=\sqrt{\tau_0}$ and $a_1=\sqrt{\tau_1}$. Fix a Riemannian metric
$\mathbf h$ on $M$. Compactness of $M$ and smoothness of the flow give
constants $c_0,C_0,R_0>0$ such that
\[
 c_0\mathbf h\leq\mathbf g_{r^2}\leq C_0\mathbf h,
 \qquad |R_{r^2}|\leq R_0,
 \qquad 0\leq r\leq a_1.
\]
Let $D=\operatorname{diam}_{\mathbf h}M$. An $\mathbf h$-minimizing
geodesic from $p$ to $q$, parametrized with constant speed on $[0,a]$,
where $a=\sqrt\tau$, has action at most
\[
 \frac{C_0d_{\mathbf h}(p,q)^2}{2a}+\frac{2R_0a^3}{3}
 \leq\frac{C_0D^2}{2a_0}+\frac{2R_0a_1^3}{3}=:C_1.
\]
Lemma~\ref{lem:per-reparametrizacion-minimizantes-L} provides a minimizer
$\alpha(r)=\gamma(r^2)$. With $\mathbf A=\alpha'$ and
$E(r)=|\mathbf A(r)|_{\mathbf g_{r^2}}^2$, formula
\eqref{eq:per-L-reparametrizada} gives
\begin{equation}
 \int_0^a E(r)\,dr
 \leq 2C_1+\frac{4R_0a_1^3}{3}=:C_2.
 \label{eq:per-energia-minimizantes-compacto}
\end{equation}
Thus some $r_*\in[a/2,3a/4]$ satisfies $E(r_*)\leq4C_2/a_0$.
Equation~\eqref{eq:per-geodesica-L-regular} and the identity
$\partial_r\mathbf g_{r^2}=4r\operatorname{Ric}_{r^2}$ yield
\[
 E'(r)=4r^2\langle\nabla^{\mathbf g_{r^2}}R_{r^2},\mathbf A\rangle_{\mathbf g_{r^2}}
       -4r\operatorname{Ric}_{r^2}(\mathbf A,\mathbf A).
\]
The uniform bounds on $\nabla R$ and $\operatorname{Ric}$ give a constant
$C_3>0$, independent of the minimizer, such that
\[
 |E'(r)|\leq C_3(1+E(r)),\qquad
 \left|\frac{d}{dr}\log(1+E(r))\right|\leq C_3.
\]
Integrating between $r_*$ and any $r\in[0,a]$, we obtain
\begin{equation}
 1+E(r)\leq\left(1+\frac{4C_2}{a_0}\right)e^{C_3a_1}.
 \label{eq:per-cota-velocidad-minimizantes}
\end{equation}
In particular, $|\mathbf A|_{\mathbf h}$ is uniformly bounded, and all
initial data $v=\mathbf A(0)/2$ lie in a fixed compact closed ball of $T_pM$.

We construct upper supports using comparison curves. Choose a finite
atlas $(U_i,\varphi_i)$ and open sets $V_i\Subset U_i$ covering $M$.
By \eqref{eq:per-cota-velocidad-minimizantes}, there is
$d\in(0,a_0/4)$, independent of the minimizer and the chart, such that,
if $q\in V_i$, the segment $\alpha([a-d,a])$ lies in a fixed compact
subset of $U_i$. Indeed, its $\mathbf h$-length is at most
$d\sup|\mathbf A|_{\mathbf h}$; it suffices to make this quantity smaller
than the distance from $\overline{V_i}$ to $M\setminus U_i$ for every
chart. The coordinate images of these segments remain a positive
distance from $\mathbb R^m\setminus\varphi_i(U_i)$.

Fix $z=(q,\tau)$, a minimizer $\alpha$ ending at $z$, and one such chart.
Write $q_i=\varphi_i(q)$ and define
\[
 b=a-d,\qquad \xi(t)=\varphi_i(\alpha(b+dt)),\qquad 0\leq t\leq1.
\]
For $y\in\varphi_i(U_i)$ and $0<\sigma\leq\tau_{\max}$ near $q_i$ and $\tau$,
respectively, put
\[
 \eta=y-q_i,\qquad c=\sqrt\sigma,\qquad h=c-b,
 \qquad x(t)=\xi(t)+t\eta,\qquad W(t)=\xi'(t)+\eta.
\]
Shrink these neighborhoods so that $d/2<h<3d/2$ and
$x([0,1])\subset\varphi_i(U_i)$. Define
\[
 \alpha_{y,\sigma}(r)=
 \begin{cases}
  \alpha(r),&0\leq r\leq b,\\
  \varphi_i^{-1}\bigl(x((r-b)/h)\bigr),&b\leq r\leq c.
 \end{cases}
\]
The two segments agree at $r=b$. This is an $H^1$ curve from $p$ to
$\varphi_i^{-1}(y)$, and $s\mapsto\alpha_{y,\sigma}(\sqrt s)$ belongs to
$\mathcal A_{p,\varphi_i^{-1}(y)}(\sigma)$.

Let $G(r,x)$ be the matrix of $\mathbf g_{r^2}$ in this chart, and put
\[
 P(r,x)=2r^2R_{r^2}(\varphi_i^{-1}(x)).
\]
The action of the comparison curve is
\begin{equation}
 \begin{split}
 \Phi_z(y,\sigma)
 :={}&\int_0^b\left(\frac12|\mathbf A(r)|_{\mathbf g_{r^2}}^2
                  +2r^2R_{r^2}(\alpha(r))\right)dr\\
 &+\int_0^1\left(
       \frac{G(b+ht,x(t))[W(t),W(t)]}{2h}
       +hP(b+ht,x(t))\right)dt.
 \end{split}
 \label{eq:per-apoyo-comparacion-explicito}
\end{equation}
By the definition of $L$,
\[
 L(\varphi_i^{-1}(y),\sigma)\leq\Phi_z(y,\sigma),
 \qquad \Phi_z(q_i,\tau)=L(q,\tau).
\]

The function $\Phi_z$ is smooth in $(y,\sigma)$. For $e\in\mathbb R^m$,
its spatial derivatives are
\begin{align*}
 D_y\Phi_z[e]
 &=\int_0^1\left(
  \frac{t}{2h}D_xG[e](W,W)+\frac1hG(W,e)+htD_xP[e]
  \right)dt,\\
 D_y^2\Phi_z[e,e]
 &=\int_0^1\left(
  \frac{t^2}{2h}D_x^2G[e,e](W,W)
  +\frac{2t}{h}D_xG[e](W,e)+\frac1hG(e,e)
  +ht^2D_x^2P[e,e]\right)dt.
\end{align*}
All coefficients in these formulas are evaluated at $(b+ht,x(t))$.
Since $b$ and $\xi$ are fixed when differentiating the support associated
with $z$,
\[
 \partial_\sigma\Phi_z
 =\frac1{2\sqrt\sigma}\int_0^1\left(
  -\frac{G(W,W)}{2h^2}+\frac{t}{2h}\partial_rG(W,W)
  +P+ht\partial_rP\right)dt.
\]
The coefficient bounds on compact sets, the estimate $|\xi'|\leq Cd$,
and the inequalities $d/2<h<3d/2$ and $\sqrt\sigma>a_0/2$ give constants
$A,H>0$ such that
\begin{equation}
 |D_{(y,\sigma)}\Phi_z|\leq A,
 \qquad D_y^2\Phi_z[e,e]\leq H|e|^2.
 \label{eq:per-cotas-derivadas-apoyos}
\end{equation}
These constants and the radii of definition are uniform over the finite
atlas and independent of the minimizer. In particular, each point has a
convex neighborhood $\mathcal W$ in space--time coordinates such that
all supports constructed in the same chart from $z\in\mathcal W$
are defined on
$\mathcal W$ and satisfy \eqref{eq:per-cotas-derivadas-apoyos}.

For $z,z'\in\mathcal W$, the support properties imply
\[
 L(z')-L(z)\leq\Phi_z(z')-\Phi_z(z)\leq A|z'-z|.
\]
Interchanging $z$ and $z'$ gives
$|L(z')-L(z)|\leq A|z'-z|$. Thus $L$ is locally Lipschitz.
On each time slice, the second bound in
\eqref{eq:per-cotas-derivadas-apoyos} and
Lemma~\ref{lem:per-apoyos-semiconcavidad} show that $L(\cdot,\tau)$ is
locally semiconcave. Since $\tau\geq\tau_0>0$ and $L$ is bounded on
compact sets, the identity $l=L/(2\sqrt\tau)$ gives both properties for $l$.

Define the action of the solution with initial datum $v$ by
\[
 \mathscr A(v,\sigma)
 :=\int_0^{\sqrt\sigma}\left(
   \frac12|\alpha_v'(r)|_{\mathbf g_{r^2}}^2
   +2r^2R_{r^2}(\alpha_v(r))\right)dr.
\]
Lemma~\ref{lem:per-exponencial-L} implies that this function is smooth
for $\sigma>0$. If $(q_j,\tau_j)\to(q,\tau)$ and $v_j$ is the initial
datum of a minimizer ending at $(q_j,\tau_j)$, estimate
\eqref{eq:per-cota-velocidad-minimizantes} yields a subsequence with
$v_j\to v$. Smooth dependence and continuity of $L$ give
\[
 \mathcal L\exp_p^\tau(v)=q,
 \qquad
 \mathscr A(v,\tau)
 =\lim_{j\to\infty}\mathscr A(v_j,\tau_j)
 =\lim_{j\to\infty}L(q_j,\tau_j)=L(q,\tau).
\]
Hence $v$ determines a minimizer ending at $(q,\tau)$.

Suppose this minimizer is unique, with initial datum $v$, and
$d(\mathcal L\exp_p^\tau)_v$ is invertible. The inverse function theorem,
applied to
\[
 (w,\sigma)\longmapsto(\mathcal L\exp_p^\sigma(w),\sigma),
\]
gives a smooth inverse between neighborhoods of $(v,\tau)$ and
$(q,\tau)$. Every minimizer ending sufficiently near $(q,\tau)$ has its
initial datum in the domain of this inverse. Otherwise, a sequence of
initial data outside that neighborhood would have, by the preceding
paragraph, a subsequence converging to the initial datum of a minimizer
ending at $(q,\tau)$. Uniqueness would force the limit to equal $v$, a
contradiction. The local inverse and existence of minimizers therefore
give uniqueness throughout this neighborhood. After shrinking the
domains, the spatial differential remains invertible. Both conditions
of Definition~\ref{def:per-lugar-corte-L} then hold, and
\[
 \mathcal U=\{(q,\tau)\in M\times(0,\tau_{\max}]
                  \mid q\in U_p(\tau)\}
\]
is open in the space--time domain. On this set, $L$ is the composition
of $\mathscr A$ with the inverse above, so $L$ and $l$ are smooth.

Finally, fix $\tau\in(0,\tau_{\max}]$. If $L(\cdot,\tau)$ is
differentiable at $q$, each minimizer ending at $q$ supplies a smooth
upper support $\Phi_z(\cdot,\tau)$ agreeing with $L$ at $q$. Their
differentials at $q$ therefore agree. The first variation formula
\eqref{eq:per-primera-variacion-L-completa}, applied to the final segment,
whose initial endpoint is fixed, yields, for every $Y\in T_qM$,
\[
 d(L(\cdot,\tau))_q(Y)
 =2\sqrt\tau\,\langle\dot\gamma(\tau),Y\rangle_{\mathbf g_\tau}.
\]
All minimizers have the same final velocity. Backward uniqueness for
the regular system \eqref{eq:per-geodesica-L-regular} implies that they
coincide. If $q$ is also not a critical value of
$\mathcal L\exp_p^\tau$, the unique minimizer has regular initial
datum; the preceding argument gives $q\in U_p(\tau)$. Consequently,
\[
 \operatorname{Cut}_{\mathcal L}(p,\tau)
 \subseteq
 \operatorname{CritVal}(\mathcal L\exp_p^\tau)
 \cup\{q\in M\mid L(\cdot,\tau)\text{ is not differentiable at }q\}.
\]
The first set is null by Sard's theorem and the second by Rademacher's
theorem. Riemannian measure has a smooth positive density in each chart,
so both sets are also $\lambda_{\mathbf g_\tau}$--null. This proves the
claim about the cut locus.
\end{proof}

\begin{lemma}[Semiconcavity and the distributional Hessian]
\label{lem:per-semiconcavidad-hessiana-distribucional}
Let $(N,\mathbf{h})$ be a smooth Riemannian manifold of dimension $m$ without boundary, and let $u$
be locally semiconcave. Then $u\in W^{1,\infty}_{\mathrm{loc}}(N)$, its
distributional Hessian is a symmetric tensor with values in Radon
measures, and the singular part of $\operatorname{Hess}_{\mathbf{h}}u$ is negative
semidefinite. In particular, the singular part of
$\Delta_{\mathbf{h}}u$ is a nonpositive measure.
\end{lemma}

\begin{proof}
Fix a chart and write $\widehat u=u\circ\varphi^{-1}$.
Shrink its image to a convex open set $V$ on which
$v(x)=\widehat u(x)-C|x|^2/2$ is concave. By
Lemma~\ref{lem:per-apoyos-semiconcavidad}, $\widehat u$ is locally
Lipschitz. Theorem~\ref{rademacher euclidiano} gives its almost everywhere
derivatives, which are locally bounded. They are also its weak derivatives:
for a test function $\psi\in C_c^\infty(V)$, transfer the quotient
$[\widehat u(x+te_j)-\widehat u(x)]/t$ to $\psi$ by a change of variables.
The Lipschitz bound and dominated convergence allow us to let $t\to0$ and give
\[
 \int_V D_j\widehat u\,\psi\,d\lambda_m
 =-\int_V\widehat u\,\partial_j\psi\,d\lambda_m.
\]
This proves that $u\in W^{1,\infty}_{\mathrm{loc}}(N)$.

Let us prove that the second derivatives of $v$ are measures. For a
constant vector $e\in\mathbb R^m$, concavity implies
$v(x+te)+v(x-te)-2v(x)\leq0$ whenever the segment lies in $V$.
If $\psi\in C_c^\infty(V)$ is nonnegative and $t\ne0$ is sufficiently
small, integrate this inequality after multiplying by $\psi/t^2$.
Transferring the increments to the test function gives
\[
 \int_V v(x)
 \frac{\psi(x+te)+\psi(x-te)-2\psi(x)}{t^2}\,d\lambda_m(x)
 \leq0.
\]
The quotient converges uniformly to $D_{ee}\psi$, with support in a
fixed compact set. Since $v$ is locally integrable, passage to the limit
gives
\[
 \langle D_{ee}v,\psi\rangle\leq0.
\]
Thus $-D_{ee}v$ is a positive distribution.
Theorem~\ref{teo: distribuciones orden cero y medidas} represents it
by a positive Radon measure. In particular, if $K\subset V$ is
compact and $0\leq\chi\in C_c^\infty(V)$ equals one near $K$,
positivity gives the order-zero bound
\[
 |\langle D_{ee}v,\psi\rangle|
 \leq\|\psi\|_\infty\langle-D_{ee}v,\chi\rangle,
 \qquad\operatorname{supp}(\psi)\subseteq K.
\]
Here positivity is applied to
$\|\psi\|_\infty\chi+\psi$ and
$\|\psi\|_\infty\chi-\psi$.

Distributional derivatives commute. The polarization identity
\[
 2D_{ij}v=D_{e_i+e_j,e_i+e_j}v-D_{ii}v-D_{jj}v
\]
shows that all $D_{ij}v$ are local signed Radon measures,
in the sense of Definition~\ref{def:b4-espacios-medida-radon-con-signo}.
Denote them by $\mu_{ij}$. The matrix
$(\mu_{ij})$ is symmetric and, for every $e\in\mathbb R^m$,
\[
 \sum_{i,j=1}^m e_ie_j\mu_{ij}=D_{ee}v\leq0.
\]
This inequality also determines the sign for directions varying
with the point. To see this, on an open set $V_0$ with compact closure
contained in $V$, take the finite positive measure
\[
 \nu=\lambda_m\restriction_{V_0}
      +\sum_{i,j=1}^m|\mu_{ij}|\restriction_{V_0}.
\]
Theorem~\ref{teo:radon-nikodym}, applied to the positive
and negative variations, gives densities $a_{ij}=d\mu_{ij}/d\nu$ and
$r=d\lambda_m/d\nu$. For each vector $e\in\mathbb Q^m$, we have
$\displaystyle\sum_{i,j=1}^{m}e_ie_ja_{ij}\leq0$ outside a $\nu$-null set.
Since $\mathbb Q^m$ is countable, we can exclude a single null set
for all these directions. By density, the symmetric matrix
$\mathbf a=(a_{ij})$ is negative semidefinite at $\nu$-almost every point.

Set $S=\{x\in V_0\mid r(x)=0\}$. Then $\lambda_m(S)=0$ and,
on $V_0\setminus S$,
the preceding measures have densities $a_{ij}/r$ with respect to
$\lambda_m$. Thus their decomposition into absolutely continuous and
singular parts is
\[
 \mu_{ij}
 =\mathbf{1}_{V_0\setminus S}\frac{a_{ij}}r\,d\lambda_m
  +\mathbf{1}_S a_{ij}\,d\nu.
\]
The measures in this equality are restricted to $V_0$, and the first
coefficient is extended by zero on $S$. The second summand is the singular part, and its matrix has the sign already
established. Since
$D_{ij}\widehat u=\mu_{ij}+C\delta_{ij}\,d\lambda_m$, the quadratic
term does not change this singular part.

We now pass to the Riemannian Hessian, taking care with the measure used to
identify functions with distributions. In this chart,
$D_{ij}\widehat u$ is a measure written with the Euclidean convention,
whereas
$d\lambda_{\mathbf h}=\sqrt{\det(\mathbf h)}\,d\lambda_m$.
If $\mathsf H_{ij}$ denotes the coefficient of the distributional Hessian
when represented by measures, then
\[
 \mathsf H_{ij}
 =\sqrt{\det(\mathbf h)}\,D_{ij}\widehat u
  -\left(\sum_{k=1}^m\Gamma_{ij}^kD_k\widehat u\right)
    d\lambda_{\mathbf h}.
\]
Indeed, for a smooth function this is the coordinate formula
$\operatorname{Hess}_{\mathbf h}u=\nabla du$ multiplied by the volume
element. For the present function, regularize $\widehat u$ by
convolution in a smaller chart $(U,\phi)$. The approximations converge in
$W^{1,1}(\phi(U))$ on compact sets, and their second derivatives converge
distributionally to $D_{ij}\widehat u$. Continuity of
distributional differentiation and multiplication by smooth
coefficients allows passage to the limit in the formula. On overlaps,
we also pass to the limit in the tensor transformation law for the smooth Hessians,
since the coordinate changes and their coefficients are smooth.
The resulting measures therefore represent the intrinsic tensor
$\nabla du$ and are independent of the chart.

The term containing the Christoffel symbols has locally
bounded density. Moreover, $\lambda_m$ and $\lambda_{\mathbf h}$ have
the same null sets. Consequently,
\[
 \mathsf H_{ij}^{\mathrm s}
 =\sqrt{\det(\mathbf h)}\,\mathbf{1}_S a_{ij}\,d\nu.
\]
The square root of the determinant is positive, and $\mathbf a$ is negative
semidefinite; this proves the assertion about the singular Hessian.
Its contraction with the inverse metric is
\[
 (\Delta_{\mathbf h}u)^{\mathrm s}
 =\sqrt{\det(\mathbf h)}\,\mathbf{1}_S
  \left(\sum_{i,j=1}^m h^{ij}a_{ij}\right)d\nu\leq0.
\]
The sign follows because the trace of the product of a positive
definite matrix and a negative semidefinite matrix is nonpositive.

We also record the local control that we will need as time
varies. If a family of functions $u_\tau$ has uniform bounds on
$\|\widehat u_\tau\|_\infty$, its Lipschitz constant, and its
semiconcavity constant in a common chart, the total variations of the
Hessian measures are uniformly bounded on every compact subset of the
interior. Indeed, for $v_\tau=\widehat u_\tau-C|x|^2/2$ and the preceding cutoff function
$\chi$,
\[
 (-D_{ee}v_\tau)(K)
 \leq\langle-D_{ee}v_\tau,\chi\rangle
 \leq\|v_\tau\|_{L^\infty(\operatorname{supp}\chi)}
       \|D_{ee}\chi\|_{L^1(\mathbb R^m)}.
\]
Apply this bound to the finitely many directions $e_i$ and $e_i+e_j$ and
use polarization. The quadratic term and the Christoffel
terms are controlled by the bounds already stated. This control remains
uniform for a smooth family of metrics on a compact
interval, since their coefficients, inverses, and first derivatives
are bounded in the smaller chart.
\end{proof}

\begin{proposition}[Variational geometry of $\mathcal L$]
\label{prop:per-variaciones-L}
Let $M$ be a smooth closed connected manifold, and let $\mathbf{g}(t)$ be a smooth
Ricci flow on an interval containing $[t_0-\tau_{\max},t_0]$. For each
$q\in M$ and $0<\tau\leq\tau_{\max}$, there exists a minimizer of
$\mathcal L$. On the regular set,
\eqref{eq:per-primera-variacion-L} and
\eqref{eq:per-tiempo-l-K}--\eqref{eq:per-laplaciano-l-K} hold. On the entire
manifold, the reduced function satisfies, in the distributional sense,
\begin{equation}
 \partial_\tau l-\Delta_\tau l+|\nabla^{\mathbf{g}_\tau}l|_{\mathbf{g}_\tau}^2
 -R_\tau+\frac m{2\tau}\geq0.
 \label{eq:per-desigualdad-l-reducida}
\end{equation}
\end{proposition}

\begin{proof}
Existence and the formulas on the regular set were proved in the
preceding lemmas. By
Lemma~\ref{lem:per-regularidad-lugar-corte-L}, $l$ is locally Lipschitz in
spacetime, and the cut locus has measure zero on each slice.
Fubini's theorem then shows that
\eqref{eq:per-desigualdad-l-reducida-regular} holds almost everywhere on
$M\times(0,\tau_{\max}]$.

Fix a compact interval $J\subset(0,\tau_{\max})$. For almost every
$\tau\in J$, write the Laplacian measure of $l(\cdot,\tau)$ as
\[
 \Delta_\tau l=(\Delta_\tau l)_{\mathrm{ac}}\,
 d\lambda_{\mathbf{g}_\tau}+\mu_\tau^{\mathrm s}.
\]
Lemma~\ref{lem:per-semiconcavidad-hessiana-distribucional} gives
$\mu_\tau^{\mathrm s}\leq0$. On the regular open set, the density
$(\Delta_\tau l)_{\mathrm{ac}}$ agrees with the classical Laplacian.
Since the complement of that open set has volume measure zero, the
almost everywhere inequality obtained above implies
\[
 \partial_\tau l-(\Delta_\tau l)_{\mathrm{ac}}
 +|\nabla^{\mathbf{g}_\tau}l|_{\mathbf{g}_\tau}^2-R_\tau
 +\frac m{2\tau}\geq0
\]
at almost every point of each of these slices.

Let $\psi\geq0$ be a smooth function with compact support contained in
$M\times\operatorname{Int}(J)$. Integrating on a slice, the sign of
$-\mu_\tau^{\mathrm s}$ gives
\[
 \int_M\left(\partial_\tau l+
 |\nabla^{\mathbf{g}_\tau}l|_{\mathbf{g}_\tau}^2-R_\tau+
 \frac m{2\tau}\right)\psi\,d\lambda_{\mathbf{g}_\tau}
 -\langle\Delta_\tau l,\psi(\cdot,\tau)\rangle\geq0.
\]
The Laplacian pairing can be written without separating its two parts:
\[
 \langle\Delta_\tau l,\psi(\cdot,\tau)\rangle
 =-\int_M\langle\nabla^{\mathbf{g}_\tau}l,
 \nabla^{\mathbf{g}_\tau}\psi\rangle_{\mathbf{g}_\tau}
 \,d\lambda_{\mathbf{g}_\tau}.
\]
This identity is the weak definition of the Laplacian and holds because
$l(\cdot,\tau)$ is locally Lipschitz. Its right-hand side depends
measurably on $\tau$ and is integrable on $J$: the weak gradient of
$l$ is jointly measurable and uniformly bounded on the support
of $\psi$, and the metric coefficients are smooth. The same holds
for the remaining terms. The uniform control of the measures on
compact sets, proved at the end of the preceding lemma, also applies here,
since the semiconcavity and Lipschitz bounds for $l$ are uniform on
$J$. Thus time integration is justified directly by
Fubini's theorem, without presupposing measurability of a separate choice of the
singular measures. We obtain
\[
 \int_J\left\{\int_M\left(\partial_\tau l+
 |\nabla^{\mathbf{g}_\tau}l|_{\mathbf{g}_\tau}^2-R_\tau+
 \frac m{2\tau}\right)\psi\,d\lambda_{\mathbf{g}_\tau}
 +\int_M\langle\nabla^{\mathbf{g}_\tau}l,
 \nabla^{\mathbf{g}_\tau}\psi\rangle_{\mathbf{g}_\tau}
 \,d\lambda_{\mathbf{g}_\tau}\right\}d\tau\geq0.
\]
Finally, the weak time derivative also agrees with the almost
everywhere derivative of $l$. For the spacetime measure we use,
$\partial_\tau d\lambda_{\mathbf{g}_\tau}=R_\tau
 d\lambda_{\mathbf{g}_\tau}$; hence integration by parts in time
takes the form
\[
 \int_J\int_M(\partial_\tau l)\psi\,
 d\lambda_{\mathbf{g}_\tau}d\tau
 =-\int_J\int_M l(\partial_\tau\psi+R_\tau\psi)\,
 d\lambda_{\mathbf{g}_\tau}d\tau.
\]
The preceding inequality is therefore exactly
\eqref{eq:per-desigualdad-l-reducida} in the distributional sense. Since the
interval $J$ and the test function were arbitrary, this justifies the
formulation used below; it corresponds to the extension
developed in \cite[Theorem~7.13]{morgan2007ricci}.
\end{proof}

\begin{definition}[Reduced volume]
\label{def:per-volumen-reducido}
The reduced volume based at $(p,t_0)$ is
\begin{equation}
 \widetilde V(\tau)
 :=\int_M(4\pi\tau)^{-m/2}e^{-l(q,\tau)}
 \,d\lambda_{\mathbf{g}(t_0-\tau)}(q).
 \label{eq:per-volumen-reducido}
\end{equation}
\end{definition}

To identify the initial value of the reduced volume, we need a
quantitative comparison with the Euclidean model.

\begin{lemma}[Gaussian behavior for small times]
\label{lem:per-asintotica-gaussiana-l}
Fix $A>0$. There exist $C_A>0$ and $\tau_A>0$ such that, in normal
coordinates for $\mathbf{g}(t_0)$ centered at $p$,
\begin{equation}
 \left|l(\exp_p(\sqrt\tau\,z),\tau)-\frac{|z|_{\mathbf{g}(t_0)}^2}{4}\right|
 \leq C_A\tau
 \label{eq:per-asintotica-local-l}
\end{equation}
for $|z|\leq A$ and $0<\tau\leq\tau_A$. Moreover, there exist $c,C>0$ such that
\begin{equation}
 l(q,\tau)\geq c\frac{d_{\mathbf{g}(t_0)}(p,q)^2}{\tau}-C\tau
 \label{eq:per-cota-inferior-gaussiana-l}
\end{equation}
for every $q\in M$ and every sufficiently small $\tau$.
\end{lemma}

\begin{proof}
The curvature bound on a short interval provides a
uniform bound on $|\operatorname{Ric}_{\mathbf g_s}|_{\mathbf g_s}$
and $|R_s|$. Lemma~\ref{lem:per-distorsion-metrica-distancias}
compares $\mathbf g_s$ with $\mathbf g(t_0)$ by factors
$e^{\pm2Ks}$. On a bounded interval, the inequalities
$e^{-2Ks}\geq1-2Ks$ and
$e^{2Ks}\leq1+2Ke^{2K\tau_0}s$ allow us to choose $C,\tau_0>0$
such that
\[
 (1-C\tau)\mathbf{g}(t_0)\leq\mathbf{g}_s\leq(1+C\tau)\mathbf{g}(t_0),
 \qquad |R_s|\leq C,
 \quad 0\leq s\leq\tau\leq\tau_0.
\]
Shrink $\tau_0$ so that $C\tau_0\leq1/2$.
The first inequality means, at each point and for each vector,
\[
 (1-C\tau)|v|_{\mathbf g(t_0)}^2
 \leq|v|_{\mathbf g_s}^2.
\]
Thus
$|v|_{\mathbf g(t_0)}\leq\sqrt2\,|v|_{\mathbf g_s}$, with
a constant independent of $s$ and $\tau$.
Let $d_0=d_{\mathbf{g}(t_0)}$. For any admissible curve,
\[
 d_0(p,q)
 \leq\int_0^\tau|\dot\gamma(s)|_{\mathbf g(t_0)}\,ds
 \leq\sqrt2\int_0^\tau
       s^{1/4}|\dot\gamma(s)|_{\mathbf g_s}s^{-1/4}\,ds.
\]
Cauchy--Schwarz with the weights $s^{1/4}$ and $s^{-1/4}$ gives
\[
 d_0(p,q)^2
 \leq C\left(\int_0^\tau\sqrt s\,|\dot\gamma|_{\mathbf{g}_s}^2ds\right)
       \left(\int_0^\tau s^{-1/2}ds\right).
\]
Since the second integral is $2\sqrt\tau$ and the curvature term is
bounded by $C\tau^{3/2}$, we obtain
\eqref{eq:per-cota-inferior-gaussiana-l}.

For the local estimate, use the variable $r=\sqrt s$. If
$q=\exp_p(\sqrt\tau z)$, the radial curve
$\alpha(r)=\exp_p(rz)$ gives, by
\eqref{eq:per-L-reparametrizada},
\[
 L(q,\tau)\leq\frac{\sqrt\tau}{2}|z|^2+C_A\tau^{3/2}.
\]
Now let $\alpha$ be a minimizer. Since the scalar term is bounded
below and the preceding upper bound is $O_A(\sqrt\tau)$, we obtain the
sharper estimate
\begin{equation}
 \int_0^{\sqrt\tau}|\alpha'(r)|_{\mathbf{g}(t_0)}^2\,dr
 \leq C_A\sqrt\tau.
 \label{eq:per-energia-minimizante-tiempo-pequeno}
\end{equation}
By Cauchy--Schwarz, for $0\leq r\leq\sqrt\tau$,
\[
 d_{\mathbf{g}(t_0)}(p,\alpha(r))
 \leq\int_0^r|\alpha'(\rho)|_{\mathbf{g}(t_0)}\,d\rho
 \leq C_A\sqrt\tau.
\]
The minimizer therefore remains in the normal ball of radius
$C_A\sqrt\tau$. Let $\xi(r)$ be its normal-coordinate representation,
identified with $T_pM$. Then
$\xi(0)=0$ and $\xi(\sqrt\tau)=\sqrt\tau\,z$.
In these coordinates, smoothness and the vanishing of the first
derivatives of $\mathbf g(t_0)$ at $p$ give
$g_{ij}(t_0,\xi(r))-g_{ij}(t_0,0)=O_A(\tau)$.
The flow equation also gives
$g_{ij}(t_0-r^2,\xi(r))-g_{ij}(t_0,\xi(r))=O_A(\tau)$.
Thus, for a constant $c_A>0$ and sufficiently small $\tau$,
\[
 |\alpha'(r)|_{\mathbf g_{r^2}}^2
 \geq(1-c_A\tau)|\xi'(r)|_{\mathbf g(t_0)|_p}^2.
\]
By Cauchy--Schwarz,
\[
 \int_0^{\sqrt\tau}|\xi'(r)|_{\mathbf g(t_0)|_p}^2\,dr
 \geq\frac{|\xi(\sqrt\tau)-\xi(0)|_{\mathbf g(t_0)|_p}^2}
              {\sqrt\tau}
 =\sqrt\tau\,|z|_{\mathbf g(t_0)|_p}^2.
\]
Since $|z|\leq A$, we conclude that
\[
 \frac12\int_0^{\sqrt\tau}|\alpha'(r)|_{\mathbf g_{r^2}}^2\,dr
 \geq\frac{\sqrt\tau}{2}|z|^2-C_A\tau^{3/2}.
\]
The curvature term has absolute value at most
$C\tau^{3/2}$. Together with the radial upper bound, this gives
\[
 \left|L(\exp_p(\sqrt\tau z),\tau)
 -\frac{\sqrt\tau}{2}|z|^2\right|\leq C_A\tau^{3/2},
\]
and division by $2\sqrt\tau$ proves
\eqref{eq:per-asintotica-local-l}.
\end{proof}

\begin{theorem}[Monotonicity of the reduced volume]
\label{teo:per-monotonia-volumen-reducido}
Let $M$ be a smooth closed connected manifold, and let $\mathbf{g}(t)$ be a smooth
Ricci flow on an interval containing $[t_0-\tau_{\max},t_0]$.
The reduced volume is locally Lipschitz and nonincreasing on
$(0,\tau_{\max}]$. In particular,
\[
 \frac d{d\tau}\widetilde V(\tau)\leq0
 \quad\text{for almost every }\tau\in(0,\tau_{\max}),
 \qquad \lim_{\tau\to 0^{+}}\widetilde V(\tau)=1.
\]
Moreover, $\widetilde V$ is invariant under parabolic rescalings.
\end{theorem}

\begin{proof}
Set
$v(q,\tau)=(4\pi\tau)^{-m/2}e^{-l(q,\tau)}$. Let us explain the passage to
$e^{-l}$ in the presence of the singular part of the Hessian. For almost every
$\tau$, write
\[
 \Delta_\tau l=(\Delta_\tau l)_{\mathrm{ac}}\,d\lambda_{\mathbf{g}_\tau}
 +\mu_\tau^{\mathrm s},
 \qquad \mu_\tau^{\mathrm s}\leq0.
\]
The weak chain rule gives
$\nabla(e^{-l})=-e^{-l}\nabla l$. For a smooth test function
$\varphi$ with compact support in a chart, apply the definition of
the measure $\Delta_\tau l$ to $\varphi e^{-l}$:
\[
 \int_M\varphi e^{-l}\,d(\Delta_\tau l)
 =-\int_Me^{-l}\langle\nabla\varphi,\nabla l\rangle
       \,d\lambda_{\mathbf{g}_\tau}
  +\int_M\varphi e^{-l}|\nabla l|^2\,d\lambda_{\mathbf{g}_\tau}.
\]
The use of a Lipschitz test function here is justified by approximating it
uniformly and in $W^{1,2}(M)$ by smooth functions; uniform convergence
controls the pairing with the measure. Comparison with the definition of
$\Delta_\tau(e^{-l})$ gives
\begin{equation}
 \Delta_\tau(e^{-l})
 =e^{-l}\bigl(-(\Delta_\tau l)_{\mathrm{ac}}+|\nabla l|^2\bigr)
 d\lambda_{\mathbf{g}_\tau}-e^{-l}\mu_\tau^{\mathrm s}.
 \label{eq:per-cadena-exponencial-semiconcava}
\end{equation}
For almost every $\tau$, set
\[
 a_\tau:=\partial_\tau l-(\Delta_\tau l)_{\mathrm{ac}}
       +|\nabla l|^2-R_\tau+\frac m{2\tau}.
\]
Proposition~\ref{prop:per-variaciones-L} gives $a_\tau\geq0$
almost everywhere. The temporal chain rule gives
$\partial_\tau v=-v(\partial_\tau l+m/(2\tau))$.
By \eqref{eq:per-cadena-exponencial-semiconcava}, as spatial measures,
\[
 (\partial_\tau v+R_\tau v)\,d\lambda_{\mathbf g_\tau}
      -\Delta_\tau v
 =-v a_\tau\,d\lambda_{\mathbf g_\tau}
      +v\mu_\tau^{\mathrm s}\leq0.
\]
Both terms have the indicated sign: $v>0$, $a_\tau\geq0$, and
$\mu_\tau^{\mathrm s}\leq0$. Integration in time gives
\begin{equation}
 (\partial_\tau-\Delta_{\mathbf g_\tau}+R_\tau)v\leq0
 \label{eq:per-subsolucion-volumen-reducido}
\end{equation}
in distributions. Integration is justified as in the proof of
Proposition~\ref{prop:per-variaciones-L}: the Laplacian pairing
is expressed through weak gradients, and the remaining terms are
integrable on compact intervals bounded away from zero.

Fix $J\Subset(0,\tau_{\max})$ and
$0\leq\theta\in C_c^\infty(\operatorname{Int}J)$.
Use the space--time test function $\psi(q,\tau)=\theta(\tau)$.
Since $M$ is closed,
$\langle\Delta_\tau v,1\rangle=0$. The identity
$\partial_\tau d\lambda_{\mathbf g_\tau}
 =R_\tau d\lambda_{\mathbf g_\tau}$ gives
\[
 -\int_J\widetilde V(\tau)\theta'(\tau)\,d\tau
 =\int_J\theta(\tau)\int_M(\partial_\tau v+R_\tau v)
          \,d\lambda_{\mathbf g_\tau}\,d\tau
 \leq0.
\]
Thus the distributional derivative of $\widetilde V$ is nonpositive.
Lipschitz regularity of $l$ on $M\times J$ and compactness of $M$
imply that $\widetilde V$ is Lipschitz on $J$.
Its weak derivative agrees with its derivative almost everywhere, and,
for $\tau_1<\tau_2$ in $J$,
\[
 \widetilde V(\tau_2)-\widetilde V(\tau_1)
 =\int_{\tau_1}^{\tau_2}\widetilde V'(\tau)\,d\tau\leq0.
\]
This proves monotonicity and the assertion about the ordinary derivative.

For the initial limit, fix $A\geq1$ and use
$q=\exp_p(\sqrt\tau z)$. On each ball $|z|\leq A$,
\eqref{eq:per-asintotica-local-l} gives
$l(q,\tau)\to|z|^2/4$ uniformly. Moreover,
\[
 d\lambda_{\mathbf{g}_\tau}(q)=\tau^{m/2}J_\tau(z)\,dz,
 \qquad J_\tau(z)\longrightarrow1
\]
uniformly on that ball. Compactness of $M$ and smoothness of the initial metric give
$\operatorname{Vol}_{\mathbf{g}(t_0)}B(p,r)\leq Cr^m$ for every $r>0$:
for small radii we use normal charts, and for the remaining radii
we use the total volume. The measures of $\mathbf{g}_\tau$ and
$\mathbf{g}(t_0)$ are uniformly comparable for small $\tau$.
Divide the complement of $B(p,A\sqrt\tau)$ into the annuli
\[
 \{q\in M\mid 2^kA\sqrt\tau\leq d_0(p,q)<2^{k+1}A\sqrt\tau\},
 \qquad k\in\mathbb N_0.
\]
The bound \eqref{eq:per-cota-inferior-gaussiana-l} then gives
\[
 \int_{M\setminus B_{\mathbf{g}(t_0)}(p,A\sqrt\tau)}v\,d\lambda_{\mathbf{g}_\tau}
 \leq C\sum_{k=0}^{\infty}(2^{k+1}A)^m e^{-c4^kA^2}.
\]
This series tends to zero as $A\to\infty$, uniformly for
small $\tau$: the exponential absorbs the polynomial factor and bounds the
series by a multiple of
$\displaystyle \sum_{k=0}^{\infty}e^{-(c/2)4^kA^2}$.
First let $\tau\to 0^{+}$ with $A$ fixed, and then let
$A\to\infty$. Dominated convergence gives
\[
 \lim_{\tau\to 0^{+}}\widetilde V(\tau)
 =(4\pi)^{-m/2}\int_{\mathbb R^m}e^{-|z|^2/4}\,dz=1.
\]

Finally, for $Q>0$, let
$\widehat{\mathbf{g}}(\sigma)=Q\mathbf{g}(t_0+\sigma/Q)$. A change of
variables in \eqref{eq:per-longitud-L} shows that
$l_{\widehat{\mathbf{g}}}(q,Q\tau)=l_{\mathbf{g}}(q,\tau)$, while
$d\lambda_{\widehat{\mathbf{g}}(-Q\tau)}=Q^{m/2}
 d\lambda_{\mathbf{g}(t_0-\tau)}$. This factor cancels
$(Q\tau)^{-m/2}$ and proves invariance.
\end{proof}

\section{The noncollapsing theorem}
\label{sec:per-no-colapso}

A curvature bound alone does not suffice to obtain a smooth limit:
a sequence may lose volume at the scale being observed. The entropy
$\mu$ prevents this phenomenon during a finite time interval. The idea is to
compare two estimates for the same quantity. Its monotonicity provides
a lower bound determined by the initial metric; a test function
concentrated on a ball, on the other hand, shows that a very small volume
ratio would make $\mu$ arbitrarily negative.

\begin{definition}[Noncollapsing at scale $r$]
\label{def:per-no-colapso}
Let $M$ be a smooth manifold without boundary, and let $\mathbf{g}(t)$ be a
Ricci flow defined on an interval $I$. Whenever
$[t_0-r^2,t_0]\subseteq I$, set
\begin{equation}
 P^-(x_0,t_0,r):=
 B_{\mathbf{g}(t_0)}(x_0,r)\times[t_0-r^2,t_0].
 \label{eq:per-cilindro-parabolico-retrogrado}
\end{equation}
We say that the flow is $\kappa$--noncollapsed at scale $r$ at
$(x_0,t_0)$ if the implication
\[
 |\operatorname{Rm}_{\mathbf{g}(t)}|_{\mathbf{g}(t)}\leq r^{-2}
 \ \text{on }P^-(x_0,t_0,r)
 \quad\Longrightarrow\quad
 \operatorname{Vol}_{\mathbf{g}(t_0)}B_{\mathbf{g}(t_0)}(x_0,r)
 \geq\kappa r^m
\]
holds.
\end{definition}

We first isolate the two functional and geometric ingredients of the
proof.

\begin{lemma}[Lower bound for $\mu$ at bounded scales]
\label{lem:per-cota-inferior-mu-escalas}
Let $(M,\mathbf{h})$ be a closed Riemannian manifold of dimension $m$. For each
$A>0$, there exists $C=C(\mathbf{h},A)<\infty$ such that
\[
 \mu(\mathbf{h},\sigma)\geq-C
 \qquad(0<\sigma\leq A).
\]
\end{lemma}

\begin{proof}
To obtain a uniform bound as $\sigma$ tends to zero, we need to
retain the dimensional coefficient of the logarithmic term. Choose
$q>2$ in the range of the embedding of $H^1(M)$ into $L^q(M)$: if $m\geq3$, take
$q<2m/(m-2)$; if $m\leq2$, any $q<\infty$ suffices.
The Euclidean Gagliardo--Nirenberg inequality we need is
\[
 \|w\|_{L^q(\mathbb R^m)}
 \leq C\|dw\|_{L^2(\mathbb R^m)}^\theta
          \|w\|_{L^2(\mathbb R^m)}^{1-\theta},
 \qquad \theta=m\left(\frac12-\frac1q\right),
 \quad w\in C_c^\infty(\mathbb R^m).
\]
It follows from the Euclidean embedding $H^1(\mathbb R^m)\hookrightarrow L^q(\mathbb R^m)$ in the Sobolev
chapter by applying it to $w(\lambda\,\cdot)$: for each $\lambda>0$,
\[
 \|w\|_{L^q(\mathbb R^m)}
 \leq C\bigl(\lambda^{1-\theta}\|dw\|_{L^2(\mathbb R^m)}
                   +\lambda^{-\theta}\|w\|_{L^2(\mathbb R^m)}\bigr).
\]
If $w\ne0$, choose $\lambda=\|w\|_{L^2(\mathbb R^m)}/\|dw\|_{L^2(\mathbb R^m)}$; the case
$w=0$ is immediate. Applying this inequality to the localizations of
$v$ in a finite atlas of $(M,\mathbf{h})$ gives
\[
 \|v\|_{L^q(M)}
 \leq C\bigl(\|dv\|_{L^2(M,T^*M)}^{\theta}\|v\|_{L^2(M)}^{1-\theta}
                 +\|v\|_{L^2(M)}\bigr),
 \qquad \theta=m\left(\frac12-\frac1q\right)\in(0,1).
\]
Derivatives of the partition of unity contribute the second summand;
the finite number of charts allows us to use a common constant $C$.
With $\|v\|_{L^2(M)}=1$ and $E=\|dv\|_{L^2(M,T^*M)}^2$, the Jensen inequality
used in \eqref{eq:per-log-sobolev-coerciva-mu} implies
\begin{align*}
 \int_Mv^2\log v^2\,d\lambda_{\mathbf{h}}
 &\leq\frac{2q}{q-2}\log\|v\|_{L^q(M)}\\
 &\leq\frac{q\theta}{q-2}\log(1+E)+C
 =\frac m2\log(1+E)+C.
\end{align*}
For $0<\eta\leq1$, set $z=\eta(1+E)>0$. Then
\[
 \frac m2\log(1+E)-\eta E
 =-\frac m2\log\eta+\frac m2\log z-z+\eta
 \leq-\frac m2\log\eta+C_m,
\]
since the function $z\mapsto(m/2)\log z-z$ attains its maximum at $z=m/2$.
Thus we obtain a constant $C_0=C_0(M,\mathbf{h})$ such that
\begin{equation}
 \int_Mv^2\log v^2\,d\lambda_{\mathbf{h}}
 \leq\eta\int_M|dv|_{\mathbf{h}}^2\,d\lambda_{\mathbf{h}}
 -\frac m2\log\eta+C_0.
 \label{eq:per-log-sobolev-parametro-no-colapso}
\end{equation}
The variational expression for $\mu$ is
\begin{align*}
 \mathcal J_{\mathbf{h},\sigma}(v)
 ={}&4\sigma\int_M|dv|_{\mathbf{h}}^2\,d\lambda_{\mathbf{h}}
 +\sigma\int_MR_{\mathbf{h}}v^2\,d\lambda_{\mathbf{h}}
 -\int_Mv^2\log v^2\,d\lambda_{\mathbf{h}}\nonumber\\
 &-\frac m2\log(4\pi\sigma)-m.
\end{align*}
If $0<\sigma\leq1/2$, take $\eta=2\sigma$ in
\eqref{eq:per-log-sobolev-parametro-no-colapso}. The coefficient of the
energy becomes $2\sigma\geq0$, and the difference between
$\frac m2\log(2\sigma)$ and $\frac m2\log(4\pi\sigma)$ is constant. Since
$R_{\mathbf{h}}$ is bounded below, we obtain a bound independent
of $\sigma$. For $1/2\leq\sigma\leq A$, use $\eta=1$; then
$4\sigma-1\geq1$ and the remaining terms are uniformly bounded
below on this compact interval. Taking the infimum over
$\|v\|_{L^{2}(M)}=1$ proves the lemma.
\end{proof}

\begin{lemma}[Cutoff function on a ball with bounded curvature]
\label{lem:per-corte-bola-no-colapso}
There exist constants $c_m,C_m>0$ with the following property. If
$(N,\mathbf{h})$ (of dimension $m$) is complete and
$|\operatorname{Rm}_{\mathbf{h}}|_{\mathbf{h}}\leq r^{-2}$ on
$B_{\mathbf{h}}(x,r)$, there exists a Lipschitz function
$\zeta\colon N\to[0,1]$ such that
\[
 \zeta=1\ \text{on }B_{\mathbf{h}}(x,r/2),\qquad
 \operatorname{supp}\zeta\subseteq\overline{B_{\mathbf{h}}(x,r)},
 \qquad |d\zeta|_{\mathbf{h}}\leq\frac{C_m}{r}
\]
almost everywhere. If
$V=\operatorname{Vol}_{\mathbf{h}}B_{\mathbf{h}}(x,r)$ and
$\displaystyle Z=\int_N\zeta^2\,d\lambda_{\mathbf{h}}$, then
\begin{equation}
 c_mV\leq Z\leq V.
 \label{eq:per-corte-volumen-no-colapso}
\end{equation}
\end{lemma}

\begin{proof}
Let $\eta\colon[0,\infty)\to[0,1]$ be a Lipschitz function that equals one on
$[0,1/2]$, vanishes on $[1,\infty)$, and satisfies
$|\eta'|\leq3$. Take
$\zeta(y)=\eta(d_{\mathbf{h}}(x,y)/r)$. The support and
gradient properties follow because the distance function is $1$--Lipschitz.
Moreover,
\[
 \operatorname{Vol}_{\mathbf{h}}B_{\mathbf{h}}(x,r/2)
 \leq Z\leq V.
\]
The curvature assumption implies
$\operatorname{Ric}_{\mathbf{h}}\geq-(m-1)r^{-2}\mathbf{h}$ on the ball.
Relative Bishop--Gromov comparison between radii $r/2$ and $r$
gives
$\operatorname{Vol}B(x,r/2)\geq c_m\operatorname{Vol}B(x,r)$ after
rescaling so that $r=1$. This proves
\eqref{eq:per-corte-volumen-no-colapso}.
\end{proof}

\begin{theorem}[$\kappa$--noncollapsing at bounded scales]
\label{teo:per-no-colapso}
Let $(M,\mathbf{g}(t))$, $t\in[0,T)$, be a Ricci flow on a
smooth closed manifold of dimension $m$, with $T<\infty$. For each $\rho>0$, there exists
$\kappa=\kappa(m,\mathbf{g}(0),T,\rho)>0$ such that the flow is
$\kappa$--noncollapsed at every $(x_0,t_0)$ and every scale
\[
 0<r\leq\min\{\rho,\sqrt{t_0}\}
\]
for which the curvature bound in
Definition~\ref{def:per-no-colapso} holds.
\end{theorem}

\begin{proof}
Since a closed manifold has finitely many components, it suffices to
work on the component of $x_0$ and finally take the minimum of the
resulting constants. We therefore assume that $M$ is connected.

Apply Corollary~\ref{cor:per-monotonia-mu} on the interval
$[0,t_0]$ with
$\tau(s)=t_0+r^2-s$. Since $\tau(t_0)=r^2$ and
$\tau(0)=t_0+r^2\leq T+\rho^2$,
Lemma~\ref{lem:per-cota-inferior-mu-escalas} gives a constant
$\mu_*>-\infty$, independent of $(x_0,t_0,r)$, such that
\begin{equation}
 \mu(\mathbf{g}(t_0),r^2)
 \geq\mu(\mathbf{g}(0),t_0+r^2)\geq\mu_*.
 \label{eq:per-cota-inferior-mu-flujo}
\end{equation}

Set $\mathbf{h}=\mathbf{g}(t_0)$,
$B=B_{\mathbf{h}}(x_0,r)$, and $V=\operatorname{Vol}_{\mathbf{h}}B$. The
parabolic assumption implies, in particular,
$|\operatorname{Rm}_{\mathbf{h}}|\leq r^{-2}$ on $B$. Let $\zeta$ be the
function in Lemma~\ref{lem:per-corte-bola-no-colapso}, let
$\displaystyle Z=\int_M\zeta^2d\lambda_{\mathbf{h}}$, and define
$v=Z^{-1/2}\zeta$. Then $v\in H^1(M)$, $v\geq0$, and
$\|v\|_{L^2(M)}=1$. The extension of the variational formula to the sphere in $H^1(M)$,
justified in Proposition~\ref{prop:per-minimizador-mu}, allows us to use
this test function.

We now estimate each term of
$\mathcal J_{\mathbf{h},r^2}(v)$. By
\eqref{eq:per-corte-volumen-no-colapso},
\[
 4r^2\int_M|dv|_{\mathbf{h}}^2d\lambda_{\mathbf{h}}
 =\frac{4r^2}{Z}\int_M|d\zeta|_{\mathbf{h}}^2d\lambda_{\mathbf{h}}
 \leq C_m\frac VZ\leq C_m.
\]
Moreover, $|R_{\mathbf{h}}|\leq C_mr^{-2}$ on the support of $v$, so
\[
 \left|r^2\int_MR_{\mathbf{h}}v^2d\lambda_{\mathbf{h}}\right|\leq C_m.
\]
For the entropy, use $v^2=\zeta^2/Z$:
\begin{align*}
 -\int_Mv^2\log v^2d\lambda_{\mathbf{h}}
 &=\log Z-\frac1Z\int_M\zeta^2\log\zeta^2d\lambda_{\mathbf{h}}\\
 &\leq\log V+C_m,
\end{align*}
since $0\leq-\zeta^2\log\zeta^2\leq e^{-1}$ and
$Z\geq c_mV$. Thus
\begin{equation}
 \mu(\mathbf{g}(t_0),r^2)
 \leq\mathcal J_{\mathbf{h},r^2}(v)
 \leq\log\frac{V}{r^m}+C_m.
 \label{eq:per-cota-superior-mu-volumen}
\end{equation}
Combining
\eqref{eq:per-cota-inferior-mu-flujo} and
\eqref{eq:per-cota-superior-mu-volumen} gives
\[
 V\geq e^{\mu_*-C_m}r^m.
\]
We may take $\kappa=e^{\mu_*-C_m}$, decreased if necessary so that
$0<\kappa\leq1$.
\end{proof}

\begin{corollary}[Nondegenerate injectivity radius]
\label{cor:per-inyectividad-no-colapso}
Let $(M,\mathbf{g}(t))$, $t\in[0,T)$, be a Ricci flow on a
smooth closed manifold of dimension $m$. If it is $\kappa$--noncollapsed at scale $r$ at
$(x_0,t_0)$ and
$|\operatorname{Rm}_{\mathbf{g}(t)}|_{\mathbf{g}(t)}\leq r^{-2}$ on
$P^-(x_0,t_0,r)$, then there exists $c=c(m,\kappa)>0$ such that
\[
 \operatorname{inj}_{\mathbf{g}(t_0)}(x_0)\geq cr.
\]
\end{corollary}

\begin{proof}
On the final slice, we have
$|\operatorname{Rm}|\leq r^{-2}$ on $B(x_0,r)$ and
$\operatorname{Vol}B(x_0,r)\geq\kappa r^m$.
Lemma~\ref{lem:per-volumen-curvatura-inyectividad} directly gives the bound
$\operatorname{inj}(x_0)\geq c(m,\kappa)r$.
\end{proof}

Together with Shi's estimates, this corollary provides the two
families of bounds needed to apply Hamilton compactness to
curvature rescalings: control of all curvature jets and
a lower bound on the injectivity radius at the basepoint.

\section{Hamilton--Ivey pinching and \texorpdfstring{$\kappa$}{kappa}--solutions}
\label{sec:per-hamilton-ivey-kappa}

In dimension three, the curvature operator may initially have a
negative part. Hamilton--Ivey pinching shows that this part becomes
negligible relative to positive curvature as the curvature
scale tends to infinity. This is why the limiting models
of a three-dimensional singularity have nonnegative curvature
operator.

Denote the smallest eigenvalue of
$\mathcal R_{\mathbf{g}(t)}\colon\Lambda^2T_xM\to\Lambda^2T_xM$ by $\nu(x,t)$, with the
normalization $\operatorname{tr}_{\Lambda^2}\mathcal R=R/2$.

\begin{theorem}[Hamilton--Ivey pinching]
\label{teo:per-hamilton-ivey}
Let $(M^3,\mathbf{g}(t))$, $0\leq t<T<\infty$, be a Ricci flow on a
smooth closed manifold. There exist constants $A,B>0$, determined by
$\mathbf{g}(0)$, such that at every point where
$\nu<0$ and $-\nu\geq B/(1+Bt)$,
\begin{equation}
 R\geq2(-\nu)\left[
 \log\!\left(\frac{-\nu}{B}\right)+\log(1+Bt)-A\right].
 \label{eq:per-hamilton-ivey-log}
\end{equation}
In particular, for each $\varepsilon>0$, there exists
$C_\varepsilon=C_\varepsilon(\mathbf{g}(0),T)$ such that
\begin{equation}
 \nu\geq-\varepsilon|\operatorname{Rm}_{\mathbf{g}(t)}|_{\mathbf{g}(t)}
 -C_\varepsilon
 \quad\text{on }M\times[0,T).
 \label{eq:per-hamilton-ivey-epsilon}
\end{equation}
\end{theorem}

The logarithmic estimate is the Hamilton--Ivey theorem, in the
normalization $\operatorname{tr}\mathcal R=R/2$ of
\cite[Theorem~4.26]{morgan2007ricci}. One may take $A=3$ and
$\displaystyle B\geq\displaystyle\max\{1,\|\displaystyle\max\{-\nu(\cdot,0),0\}\|_\infty\}$.
Indeed, replacing $\mathbf{g}$ by $B\mathbf{g}$ and $t$ by $Bt$
makes the smallest initial eigenvalue at least $-1$; undoing this normalization
gives \eqref{eq:per-hamilton-ivey-log}. The maximum principle is
applied through spectral support functions, so it does not require
$\nu$ to be differentiable at points of multiplicity.

Let us derive the consequence we will use in the limits.
Fix $\varepsilon>0$ and a dimensional constant $c_3$ such that
$|R|\leq c_3|\operatorname{Rm}|$. If
$a=-\nu>B\exp(A+c_3/(2\varepsilon))$, the threshold in the
theorem is satisfied, and $\log(1+Bt)\geq0$ implies
\[
 \frac{c_3}{\varepsilon}a
 \leq2a\bigl(\log(a/B)-A\bigr)
 \leq R\leq c_3|\operatorname{Rm}|.
\]
Thus $a\leq\varepsilon|\operatorname{Rm}|$. At the remaining points,
the negative part of $\nu$ is bounded by
$C_\varepsilon=B\exp(A+c_3/(2\varepsilon))$.
Both alternatives give \eqref{eq:per-hamilton-ivey-epsilon}.

\begin{definition}[Three-dimensional $\kappa$--solution]
\label{def:per-kappa-solucion}
A $\kappa$--solution is a Ricci flow
$(N^3,\mathbf{h}(s))$, $s\in(-\infty,0]$, on a smooth,
connected manifold without boundary, satisfying the following properties:
\begin{enumerate}
 \item each $(N,\mathbf{h}(s))$ is complete, and the curvature is
 uniformly bounded on $N\times[a,b]$ for each
 $-\infty<a\leq b\leq0$;
 \item the flow is nonflat and has nonnegative curvature operator;
 \item the flow is $\kappa$--noncollapsed at all scales.
\end{enumerate}
The last condition means that, for every $(x,s)$ and every $r>0$ such
that the corresponding parabolic cylinder is contained in the time
interval and
$|\operatorname{Rm}|\leq r^{-2}$ on it, we have
$\operatorname{Vol}_{\mathbf{h}(s)}B_{\mathbf{h}(s)}(x,r)\geq\kappa r^3$.
\end{definition}

Noncollapsing is preserved under smooth limits. We spell out this step
because it will be needed to verify that the limit of the rescalings
is still a $\kappa$--solution.

\begin{lemma}[Stability of noncollapsing under smooth convergence]
\label{lem:per-estabilidad-no-colapso-limite}
Let
$(M_i^m,\mathbf{g}_i(t),p_i)\to
(M_\infty^m,\mathbf{g}_\infty(t),p_\infty)$ be pointed smooth convergence
on compact spacetime sets. Suppose there exist numbers
$R_i\to\infty$ such that each flow $\mathbf{g}_i$ is
$\kappa$--noncollapsed at all its points and at every scale
$0<r\leq R_i$ whose backward parabolic cylinder is contained in its
time interval. Then the limit flow is $\kappa$--noncollapsed at
every finite scale for which the corresponding cylinder is defined.
\end{lemma}

\begin{proof}
Fix $x_\infty\in M_\infty$, a time $t_0$, and $r>0$ for
which $P^-(x_\infty,t_0,r)$ is defined and
$|\operatorname{Rm}_{\mathbf{g}_\infty}|\leq r^{-2}$ holds on it.
Let $0<\rho<\sigma<r$. By completeness, the closed ball of radius
$\sigma$ in the final slice is compact. Take a convergence
embedding $\Phi_i$ defined on a neighborhood of this ball and on a larger open set
containing the curves considered below; write
$x_i=\Phi_i(x_\infty)$. The ball exhaustion property ensures
that its image contains the pointed balls of fixed radius that we need.

For large $i$, we have
\begin{equation}
 B_{\mathbf{g}_i(t_0)}(x_i,\rho)
 \subseteq\Phi_i\bigl(B_{\mathbf{g}_\infty(t_0)}(x_\infty,\sigma)\bigr).
 \label{eq:per-inclusion-bolas-limite-no-colapso}
\end{equation}
Indeed, let $\gamma_i$ be a minimizing geodesic of length less
than $\rho$ starting at $x_i$. It can be lifted through $\Phi_i$
as long as it remains in the image of the ball of radius $\sigma$.
If the lifted curve first exited that ball, its length
in the limit metric up to that instant would be at least $\sigma$.
Uniform convergence of the metrics on the closed ball implies that
the corresponding length of $\gamma_i$ is at least
$(1-\varepsilon_i)\sigma$, with $\varepsilon_i\to0$. For large $i$,
this is greater than $\rho$, a contradiction. This proves the inclusion.

Convergence of the curvatures on the compact set
$\overline{B_{\mathbf{g}_\infty(t_0)}(x_\infty,\sigma)}
\times[t_0-\rho^2,t_0]$ and the strict margin $r^{-2}<\rho^{-2}$ give
\[
 |\operatorname{Rm}_{\mathbf{g}_i}|\leq\rho^{-2}
 \quad\text{on }P^-(x_i,t_0,\rho)
\]
for large $i$, using \eqref{eq:per-inclusion-bolas-limite-no-colapso}.
Moreover, $\rho<R_i$ for all these indices. Noncollapsing of the
approximating flows yields
\begin{equation}
 \operatorname{Vol}_{\mathbf{g}_i(t_0)}
 B_{\mathbf{g}_i(t_0)}(x_i,\rho)\geq\kappa\rho^m.
 \label{eq:per-volumen-no-colapso-aproximantes}
\end{equation}
By the ball inclusion and convergence of the volume densities
on a fixed domain,
\begin{align*}
 \kappa\rho^m
 &\leq\liminf_{i\to\infty}
   \operatorname{Vol}_{\mathbf{g}_i(t_0)}B_{\mathbf{g}_i(t_0)}(x_i,\rho)\\
 &\leq\operatorname{Vol}_{\mathbf{g}_\infty(t_0)}
                  B_{\mathbf{g}_\infty(t_0)}(x_\infty,\sigma).
\end{align*}
First let $\rho\to\sigma^{-}$ and then $\sigma\to r^{-}$.
Continuity of measure from below gives
$\operatorname{Vol}_{\mathbf{g}_\infty(t_0)}B(x_\infty,r)\geq\kappa r^m$.
\end{proof}

\begin{theorem}[Rescaling limit at a singular time]
\label{teo:per-blowup-kappa-solucion}
Let $(M^3,\mathbf{g}(t))$, $0\leq t<T<\infty$, be a maximal solution on
a smooth closed manifold. Let $(x_i,t_i)$ be spacetime maxima
as in \eqref{eq:per-maximos-espaciotemporales}, with an increasing sequence of times $t_i\to T^{-}$, and
let
$Q_i=|\operatorname{Rm}_{\mathbf{g}(t_i)}|_{\mathbf{g}(t_i)}(x_i)$. Then
$Q_i\to\infty$, and a subsequence of the rescaled flows
\[
 \mathbf{g}_i(s):=Q_i\mathbf{g}(t_i+s/Q_i),
 \qquad -Q_it_i\leq s\leq0,
\],
pointed at $x_i$, converges smoothly on compact subsets of
$M_\infty\times(-\infty,0]$ to a $\kappa$--solution
$(M_\infty,\mathbf{g}_\infty(s),x_\infty)$, for some $\kappa>0$.
\end{theorem}

\begin{proof}
The continuation criterion in the previous chapter implies that
$Q_i\to\infty$; since $t_i\to T>0$, we also have $Q_it_i\to\infty$. By the
choice of the spacetime maxima,
\begin{equation}
 |\operatorname{Rm}_{\mathbf{g}_i(s)}|_{\mathbf{g}_i(s)}\leq1
 \quad(-Q_it_i\leq s\leq0),
 \qquad
 |\operatorname{Rm}_{\mathbf{g}_i(0)}|_{\mathbf{g}_i(0)}(x_i)=1.
 \label{eq:per-curvatura-reescalada-maxima}
\end{equation}

Fix $A>0$. For large $i$, the flow $\mathbf{g}_i$ is defined on
$[-A-1,0]$. Shi's estimates applied on this interval provide
uniform bounds for all curvature derivatives on $[-A,0]$.
Now take unit scale at $(x_i,0)$. In the original variables,
this corresponds to radius $Q_i^{-1/2}$. For large $i$, this radius is smaller than
$\rho$ and $\sqrt{t_i}$, and
\eqref{eq:per-curvatura-reescalada-maxima} satisfies the parabolic assumption
of Theorem~\ref{teo:per-no-colapso}. Invariance of the volume ratio
under rescaling gives
\[
 \operatorname{Vol}_{\mathbf{g}_i(0)}B_{\mathbf{g}_i(0)}(x_i,1)\geq\kappa
\]
with a common constant $\kappa>0$.
Corollary~\ref{cor:per-inyectividad-no-colapso} yields
$\operatorname{inj}_{\mathbf{g}_i(0)}(x_i)\geq c(3,\kappa)$.

The first part of the proof of
Theorem~\ref{teo:per-compactacion-hamilton}, applied to the metrics
$\mathbf{g}_i(0)$, provides a single limit manifold $M_\infty$ and
spatial embeddings on an exhaustion. Fix these embeddings.
For each $A\in\mathbb N$, the estimates on $[-A-1,0]$ allow us to
apply the second part of that proof on $[-A,0]$. Choose
nested subsequences as $A$ increases and take their diagonal. The time-dependent limits
agree on overlaps because they are computed using the
same embeddings and on the same manifold. This gives a complete ancient flow
$(M_\infty,\mathbf{g}_\infty(s),x_\infty)$ for every $s\leq0$. The unit curvature value at the basepoint passes to the limit,
so the flow is nonflat. The bound in
\eqref{eq:per-curvatura-reescalada-maxima} also passes to the limit and shows
that each slice has bounded curvature.

Let us verify the two remaining properties. Fix, for example, the
original range $\rho=1$ in Theorem~\ref{teo:per-no-colapso}. On the rescaled
interval $[-A,0]$, a scale $r$ corresponds to radius
$rQ_i^{-1/2}$ and original time
$t_i+s/Q_i$. Thus the theorem applies simultaneously at all
points of this interval whenever
\[
 0<r\leq R_{i,A}:=
 \min\left\{\sqrt{Q_i},\sqrt{Q_it_i-A}\right\}.
\]
Since $Q_i\to\infty$ and $Q_it_i\to\infty$, we have
$R_{i,A}\to\infty$. The same constant $\kappa$ works for all
rescaled flows and every bounded radius. Applying
Lemma~\ref{lem:per-estabilidad-no-colapso-limite} on each interval
$[-A,0]$ and then letting $A\to\infty$, we conclude that the limit is
$\kappa$--noncollapsed at all scales.

Finally, the eigenvalues of the curvature operator are multiplied by
$Q_i^{-1}$. Inequality
\eqref{eq:per-hamilton-ivey-epsilon} becomes
\[
 \nu_{\mathbf{g}_i(s)}
 \geq-\varepsilon|\operatorname{Rm}_{\mathbf{g}_i(s)}|_{\mathbf{g}_i(s)}
 -\frac{C_\varepsilon}{Q_i}.
\]
Letting $i\to\infty$ gives
$\nu_{\mathbf{g}_\infty(s)}\geq-\varepsilon
|\operatorname{Rm}_{\mathbf{g}_\infty(s)}|$. Since $\varepsilon>0$ is
arbitrary, the curvature operator of the limit is nonnegative. Thus
all the conditions in
Definition~\ref{def:per-kappa-solucion} hold.
\end{proof}

\section{Structural theorems for surgery and the topological conclusion}
\label{sec:per-cirugia-poincare}

The preceding results describe the limits that arise when a region
of high curvature is magnified. To continue the flow beyond a
singularity, we must identify these regions with canonical models,
cut along nearly cylindrical necks, and attach caps whose
geometry is compatible with pinching and noncollapsing.

The full construction includes the standard solution, persistence
of caps, noncollapsing across surgeries, the canonical
neighborhood theorem, and the inductive choice of parameters. These
arguments occupy Chapters~12--17 of \cite{morgan2007ricci}.
We will use their structure and existence theorems to determine the
geometric and topological effect of each surgery. The final step will be to
reconstruct the initial manifold from the components that disappear.

The model for a neck is the shrinking round cylinder
$S^2\times\mathbb R$. If $\mathbf{g}_{S^2(1)}$ is the metric of the unit
sphere, set
\[
 \mathbf{g}_{\mathrm{cil}}(s)
 =(2-2s)\mathbf{g}_{S^2(1)}+(\mathbf{d}z)^2,
 \qquad -1\leq s\leq0.
\]
The scalar curvature of the slice $s=0$ equals one.

\begin{definition}[Necks, caps, and canonical components]
\label{def:per-cuellos-casquetes}
Fix $0<\varepsilon<1$ and
$k=\lfloor\varepsilon^{-1}\rfloor$.
Let $(M^3,\mathbf{g}(t))$ be a Ricci flow and let $(x,t_0)$ be a point.

An open set $U\subseteq M$ is an $\varepsilon$--\emph{neck of scale}
$h>0$ centered at $x$ if there exists a diffeomorphism
\[
 \Psi\colon S^2\times(-\varepsilon^{-1},\varepsilon^{-1})\longrightarrow U,
 \qquad x\in\Psi(S^2\times\{0\}),
\]
such that
\[
 \bigl\|h^{-2}\Psi^*\mathbf{g}(t_0)-\mathbf{g}_{\mathrm{cil}}(0)
 \bigr\|_{C^k(\mathbf{g}_{\mathrm{cil}}(0))}<\varepsilon.
\]
It is a \emph{strong} $\varepsilon$--\emph{neck} if $\Psi$ extends to a
spacetime identification with spatial map independent of $s$,
in which, after setting $s=h^{-2}(t-t_0)$, we have
\[
 \bigl\|h^{-2}\Psi^*\mathbf{g}(t_0+h^2s)
 -\mathbf{g}_{\mathrm{cil}}(s)\bigr\|_{C^k}<\varepsilon.
\]
This last norm is taken on
$S^2\times(-\varepsilon^{-1},\varepsilon^{-1})\times[-1,0]$,
using the reference metric
$\mathbf{g}_{\mathrm{cil}}(0)+(\mathbf{d}s)^2$ to compute derivatives.
The region under consideration must exist throughout this time interval.

An $(C,\varepsilon)$--\emph{cap} is an open set $U$ diffeomorphic to
an open ball $B^3$ or to $\mathbb{RP}^3\setminus\overline{B^3}$ that decomposes into a
compact core $K$ and an end that is an $\varepsilon$--neck.
If $x\in K$, we require $R(x,t_0)>0$ and, in the normalized metric
$\widehat{\mathbf{g}}=R(x,t_0)\mathbf{g}(t_0)$,
\[
 \operatorname{diam}_{\widehat{\mathbf{g}}}K\leq C,
 \qquad \operatorname{Vol}_{\widehat{\mathbf{g}}}K\geq C^{-1},
 \qquad C^{-1}\leq\frac{R(y,t_0)}{R(x,t_0)}\leq C
 \quad(y\in K).
\]

A $C$--\emph{component} is a compact component of positive
curvature that, when normalized by the scalar curvature at any of its
points, has diameter at most $C$, scalar curvature between $C^{-1}$ and $C$,
and volume between $C^{-1}$ and $C$. An
$\varepsilon$--\emph{round component} is a component that, after
this normalization, lies at distance less than $\varepsilon$ in
$C^k$ from a round spherical space form of scalar curvature one.
\end{definition}

\begin{definition}[Strong canonical neighborhood]
\label{def:per-vecindad-canonica}
Given $C<\infty$ and $\varepsilon>0$, a point $(x,t)$ of positive scalar
curvature has a strong canonical neighborhood $(C,\varepsilon)$ if
one of the following alternatives holds:
\begin{enumerate}
 \item $x$ is the center of a strong $\varepsilon$--neck whose scale is
 comparable to $R(x,t)^{-1/2}$;
 \item $x$ belongs to the core of an $(C,\varepsilon)$--cap;
 \item the component of $x$ is a $C$--component;
 \item the component of $x$ is $\varepsilon$--round.
\end{enumerate}
In all cases, the quantitative estimates
\begin{equation}
 |\nabla R|(x,t)\leq C R(x,t)^{3/2},
 \qquad |\partial_tR|(x,t)\leq C R(x,t)^2,
 \label{eq:per-cotas-vecindad-canonica}
\end{equation}
and a volume lower bound at scale $R(x,t)^{-1/2}$ of the form
$c(C,\varepsilon)R(x,t)^{-3/2}$ are included.
\end{definition}

To close off a neck, fix a rotationally symmetric metric
$\mathbf{g}_{\mathrm{std}}$ on $\mathbb R^3$ with nonnegative curvature
operator and positive scalar curvature that agrees outside a compact set
with the normalized cylinder. Its flow is the model for the caps
that are inserted.

\begin{theorem}[Standard solution]
\label{teo:per-solucion-estandar}
For the fixed standard metric $\mathbf{g}_{\mathrm{std}}$, there exists a unique
complete solution of the Ricci flow starting from it, defined on
$[0,1)$ and with uniformly bounded curvature on each
$\mathbb R^3\times[0,T_0]$, $T_0<1$. On these intervals, it is
asymptotic to all orders, away from the cap, to the round
cylinder of scalar curvature $1/(1-t)$.

There exist $r_{\mathrm{std}}>0$ and $\kappa_{\mathrm{std}}>0$ such that it is
$\kappa_{\mathrm{std}}$--noncollapsed at scales smaller than
$r_{\mathrm{std}}$ whose cylinders are contained in $[0,1)$.
For each sufficiently small $\varepsilon>0$, there exist
$C_\varepsilon<\infty$ and $R_\varepsilon<\infty$ such that the points
with $R(x,t)\geq R_\varepsilon$ have strong canonical neighborhoods
$(C_\varepsilon,\varepsilon)$. The constants are chosen after fixing
$\varepsilon$ and the standard metric.
\end{theorem}

Existence, rotational symmetry, behavior at infinity, and
uniqueness are proved in
\cite[Theorem~12.5 and Theorems~12.28 and 12.32]{morgan2007ricci}.
Uniqueness uses a harmonic map flow and is the noncompact version of the gauge mechanism
studied in Chapter~\ref{cap:flujo-ricci}.

\begin{definition}[Ricci flow with surgery]
\label{def:per-flujo-con-cirugia}
A Ricci flow with surgery, initially defined up to a time
$T_*\leq\infty$, consists of strictly increasing surgery times
$t_j$, with $j$ in an initial subset of $\mathbb N$, and of $t_0=0$.
The list may be finite, empty, or infinite. In the last case, it may,
a priori, accumulate at $T_*$. Between consecutive times there are
closed manifolds $M_j$ and smooth flows
$\mathbf{g}_j(t)$ on $M_j\times[t_j,t_{j+1})$; if there is a last surgery
time, the remaining interval ends at $T_*$. At a surgery time
$t_{j+1}$, finite disjoint collections of strong necks are chosen;
they are cut along their central spheres, the ends and
components prescribed by the procedure are removed, and the new
boundaries are closed off with caps obtained from the standard solution at the same
scale. There is an open set, called the continuing region, that is identified
isometrically before and after surgery, and the new metric agrees
there with the limit of the previous metric. The flow becomes extinct at
$T_{\mathrm{ext}}$ if all components have been discarded, and we set
$M(t)=\varnothing$ for $t\geq T_{\mathrm{ext}}$.
\end{definition}

Cutting must be performed at scales sufficiently separated
from the canonical neighborhood scale. The following parameters record this
choice and the properties that must be preserved when it is repeated.

\begin{definition}[Standard control of surgery parameters]
\label{def:per-control-estandar-cirugia}
Fix small $\varepsilon>0$ and the canonical neighborhood constant $C$.
Standard control consists of positive nonincreasing
functions
\[
 r(t),\qquad\kappa(t),\qquad\delta(t)
\]
and a surgery scale
\[
 \rho(t)=\delta(t)r(t),
 \qquad h(t)=h(\rho(t),\delta(t))\leq\delta(t)^2r(t).
\]
A flow with surgery is subject to this control if it satisfies:
\begin{enumerate}
 \item Hamilton--Ivey pinching;
 \item noncollapsing with constant $\kappa(t)$ at scales
 $0<\rho\leq\rho_*$, where $\rho_*>0$ is fixed, on cylinders that
 survive without surgery throughout their time interval; in the
 normalization of the source, one may take $\rho_*=\varepsilon$;
 \item the strong canonical neighborhood assumption at every point with
 $R(x,t)\geq r(t)^{-2}$;
 \item each surgery is performed on strong necks of accuracy
 $\delta(t)$, centered at points with $R=h(t)^{-2}$, and the post-surgery metric
 is constructed using the standard cap of scale $h(t)$;
 \item the functions $r$, $\kappa$, and $\delta$ are positive and nonincreasing;
 in particular, for each $T<\infty$, their values on $[0,T]$ are bounded
below by the positive numbers $r(T)$, $\kappa(T)$, and $\delta(T)$.
 The choice of scale $h$ also satisfies
 $\displaystyle h_T:=\inf_{t\in[0,T]}h(t)>0$.
\end{enumerate}
The full formulation includes the seven compatibility conditions
for the spacetime, the continuing and disappearing regions, and
maximality between consecutive surgery times in
\cite[Chapter~15, Sections~1 and 3]{morgan2007ricci}.
\end{definition}

The first structural result describes the ancient models that arise
at the curvature scale.

\begin{theorem}[Structure of three-dimensional $\kappa$--solutions]
\label{teo:per-estructura-kappa-soluciones}
There exists $\overline\varepsilon>0$ such that, for each $\kappa>0$ and
$0<\varepsilon\leq\overline\varepsilon$, one can choose
$C=C(\kappa,\varepsilon)$ with the following property. In an
orientable three-dimensional $\kappa$--solution, every point has, at its
curvature scale, an $\varepsilon$--neck, an
$(C,\varepsilon)$--cap, a $C$--component, or an
$\varepsilon$--round component. The approximations include the curvature, derivative, and
volume estimates of Definition~\ref{def:per-vecindad-canonica}. When the
neck alternative holds, the neck can be chosen strong. The scalar curvature is
positive at every point: if it vanished, the strong maximum principle for
$\partial_sR=\Delta R+2|\operatorname{Ric}|^2$, connectedness, and
uniqueness of the complete flow with bounded curvature would force the
flat solution, which is excluded by definition.
\end{theorem}

This is the qualitative description obtained in
\cite[Chapter~9, in particular Theorem~9.93 and
Corollaries~9.94--9.95]{morgan2007ricci}. Orientability excludes the
cylindrical quotient $\mathbb{RP}^2\times\mathbb R$.

\begin{theorem}[Canonical neighborhood theorem]
\label{teo:per-vecindades-canonicas}
Fix a smooth initial metric $\mathbf{g}_0$ on a closed,
orientable three-dimensional manifold, a time horizon $T<\infty$, and
sufficiently small $\varepsilon>0$. There exist
$r=r(\mathbf{g}_0,T,\varepsilon)>0$ and
$C=C(\mathbf{g}_0,T,\varepsilon)<\infty$ such that, in the smooth flow
starting from $\mathbf{g}_0$, every point before $T$ with
$R(x,t)\geq r^{-2}$ has a strong canonical neighborhood
$(C,\varepsilon)$. The threshold is also chosen so that the strong
necks have the required time interval available.

For the flows with surgery of
Theorem~\ref{teo:per-existencia-control-cirugia}, this property is
preserved with the time-dependent threshold $r(t)$ chosen together with the
surgery parameters. After normalization by $R(x,t)$, the high-curvature regions
are approximated by regions of $\kappa$--solutions or
by regions of the standard solution modeling a recent surgery.
\end{theorem}

The proof is by contradiction. One chooses the first points where the
conclusion fails, normalizes their curvature, and applies noncollapsing, curvature
control at bounded distance, and Hamilton compactness. If the
rescaled times extend to $-\infty$, the limit is a
$\kappa$--solution, and the description of its neighborhoods gives a contradiction. If a
surgery prevents this extension, the limit is modeled on the standard
solution, yielding the same contradiction. The complete argument,
including the curvature estimate at bounded distance and the case near
a surgery cap, is developed in
\cite[Chapters~10--12 and 17]{morgan2007ricci}.

\begin{theorem}[Existence and control of the flow with surgery]
\label{teo:per-existencia-control-cirugia}
Let $M$ be a smooth closed connected orientable three-dimensional manifold, and
let $\mathbf{g}_0$ be a smooth metric. There exist control parameters
$r(t)$, $\kappa(t)$, $\delta(t)$, and $h(t)$ and a Ricci flow with surgery
starting from $\mathbf{g}_0$, defined for every $t\geq0$ up to possible
extinction. The flow satisfies Hamilton--Ivey pinching, noncollapsing,
and Theorem~\ref{teo:per-vecindades-canonicas}; each surgery is
performed on a finite collection of strong necks of scale $h(t)$. The
surgery times do not accumulate in any compact time interval.
\end{theorem}

To apply the formulation of
\cite[Theorem~15.9 and Corollary~15.10]{morgan2007ricci} to an arbitrary
metric, first multiply $\mathbf{g}_0$ by a constant to
normalize the curvature and the volume of unit balls. Parabolic
invariance of the flow then allows this normalization to be undone.
Persistence of noncollapsing across caps is proved in
Chapter~16, and the strong canonical neighborhood assumption and absence of
accumulation of surgeries are completed in Chapter~17 of the same source.
Let us make the absence of accumulation precise, as used in
\cite[Lemma~17.12]{morgan2007ricci}. Fix $T<\infty$.
Pinching gives a bound $R\geq-C_T$, so the total
volume $V(t)$ satisfies $V'(t)\leq C_TV(t)$ during the smooth intervals.
The cap construction and the bound $h(t)\geq h_T>0$ give a
number $v_T>0$ such that each neck cut decreases the net volume
by at least $v_T$ after the cap volumes have been added.
Upon multiplication by $e^{-C_Tt}$, smooth evolution makes the weighted
volume nonincreasing, and each cut produces a decrease of at least
$e^{-C_TT}v_T$. Thus the total number $S(T)$ of cuts satisfies
\[
 S(T)\leq\frac{e^{C_TT}V(0)}{v_T}.
\]
There may also be times when components are merely discarded,
without any neck being cut. Each cut increases the number of components by at
most one. If there were initially $N_0$ components, the number $D(T)$ of
discarded components satisfies $D(T)\leq N_0+S(T)$.
Every surgery time contains a cut or a discard, so there are
at most $S(T)+D(T)$ such times in $[0,T]$.

The topological effect of each surgery is also fully controlled:
one cuts along spheres and discards only components
covered by necks and caps.

\begin{theorem}[Topological reconstruction from surgeries]
\label{teo:per-reconstruccion-topologica-cirugia}
Let $M$ be a smooth closed connected orientable three-dimensional manifold. If
a flow with surgery subject to the control of
Theorem~\ref{teo:per-existencia-control-cirugia} becomes extinct, then
$M$ is diffeomorphic to a finite connected sum of spherical space forms and
copies of $S^2\times S^1$.
\end{theorem}

\begin{proof}
Cutting along an embedded sphere and capping the two resulting
boundaries with balls gives two cases. If the sphere separates the component,
one obtains two closed manifolds whose connected sum reconstructs the
original component. If it does not separate it, the result is a closed manifold
$N$, and the original component is $N\#(S^2\times S^1)$; orientability
excludes the nonorientable sphere bundle over the circle.

The classification of components discarded by the standard
procedure gives spherical space forms, $S^2\times S^1$, and
$\mathbb{RP}^3\#\mathbb{RP}^3$. The last case is already a connected sum
of two spherical space forms. Control of these alternatives,
including the regions that disappear upon reaching the singular time,
is Proposition~15.3 and Corollary~15.4 of
\cite{morgan2007ricci}. Thus, at each step, the preceding manifold is
reconstructed from the subsequent components by adding only
summands of the indicated types.

Let $t_N$ be the extinction time. The existence and control theorem gives
only finitely many surgery times in $[0,t_N]$.
Start with the empty slice after $t_N$: the components
discarded at that step have the decomposition described above.
Suppose that all components immediately after
$t_j$, with $j\in\{1,\ldots,N\}$, have been reconstructed. Reverse the cuts at $t_j$
using the two operations in the first paragraph and add the components
discarded at that time. The result is again a finite connected sum
of spherical space forms and copies of $S^2\times S^1$.
This proves the step back to the preceding slice. Upon reaching $t_0=0$,
we obtain the decomposition of $M$.
\end{proof}

For the Poincaré conjecture, it remains to know that, under the
corresponding assumption on the fundamental group, the flow disappears in finite
time.

\begin{theorem}[Finite-time extinction]
\label{teo:per-extincion-finita-cirugia}
Let $M$ be a smooth closed connected orientable three-dimensional manifold, and
let $(M(t),\mathbf{g}(t))$ be a flow with surgery subject to the control of
Theorem~\ref{teo:per-existencia-control-cirugia}. If
\begin{equation}
 \pi_1(M)\cong G_1*\cdots*G_a*
 \underbrace{\mathbb Z*\cdots*\mathbb Z}_{b\text{ factors}},
 \qquad a,b\in\mathbb N_0,\quad |G_j|<\infty\ (j\in\{1,\ldots,a\}),
 \label{eq:per-grupo-extincion}
\end{equation}
then the flow becomes extinct in finite time.
\end{theorem}

The complete proof is in Chapters~18--19 of
\cite{morgan2007ricci}, together with the correction to its Section~19.2 in
\cite{morganTian2015correction}. First, the least area of spheres nontrivial
in $\pi_2$ is controlled through a forward difference differential
inequality. Then, when the remaining components have
$\pi_2=0$ and finite fundamental group, one uses a nontrivial class in
$\pi_3$ represented by families of loops and least disk areas. The
same differential structure forces the minimax quantity to become
negative in finite time if the component were to persist, which is
impossible. The modifications at a surgery time are compatible with
these quantities by semicontinuity.

\begin{theorem}[Poincaré conjecture]
\label{teo:per-poincare}\index{Poincare conjecture@Poincaré conjecture}
Every smooth closed simply connected three-dimensional manifold is
diffeomorphic to $S^3$.
\[
 M^3\text{ closed and simply connected}
 \quad\Longrightarrow\quad M^3\cong S^3.
\]
\end{theorem}

\begin{proof}
Let $M$ be a manifold as in the statement. It is orientable: otherwise,
its orientation cover would be a connected two-sheeted cover and would give
a nontrivial quotient of $\pi_1(M)$ isomorphic to $\mathbb Z_2$. Choose a
smooth metric $\mathbf{g}_0$ and apply
Theorem~\ref{teo:per-existencia-control-cirugia}. Since $\pi_1(M)=0$,
assumption \eqref{eq:per-grupo-extincion} holds with the empty free
product. Theorem~\ref{teo:per-extincion-finita-cirugia} shows that the
flow becomes extinct in finite time.

By Theorem~\ref{teo:per-reconstruccion-topologica-cirugia}, there exist
spherical space forms $N_j=S^3/\Gamma_j$ and an integer $b\geq0$ such
that
\[
 M\cong N_1\#\cdots\#N_a\#
 \underbrace{(S^2\times S^1)\#\cdots\#(S^2\times S^1)}_{b\text{ summands}}.
\]
The Seifert--van Kampen theorem gives
\[
 \pi_1(M)\cong\Gamma_1*\cdots*\Gamma_a*
 \underbrace{\mathbb Z*\cdots*\mathbb Z}_{b\text{ factors}}.
\]
By the normal form theorem for free products, this group can
be trivial only if $b=0$ and each $\Gamma_j$ is trivial. Consequently,
each $N_j$ is diffeomorphic to $S^3$. Since $S^3$ is the identity for the connected
sum of closed orientable three-dimensional manifolds, we conclude
that $M\cong S^3$.
\end{proof}

\part{Infinite-dimensional geometry}
\chapter*{Introduction to Part VI}
\addcontentsline{toc}{chapter}{Introduction to Part VI}
\markboth{Part VI. Infinite-dimensional geometry}{Introduction to Part VI}

A change of perspective arises naturally in many problems of geometric analysis. So far, we have studied functions, sections, metrics, and connections as objects defined on a manifold. From this part onward, we will also regard them as points of spaces with their own geometric structure. A space of maps, the set of Riemannian metrics, or the space of connections can, under suitable hypotheses, behave as an infinite-dimensional manifold.

Working rigorously with these spaces requires a careful choice of regularity. Banach spaces allow the use of versions of the inverse function, rank, and regular preimage theorems. Hilbert spaces add an inner product structure and the Riesz identification. Banach--Finsler structures provide tangent norms and distances that will be useful in studying flows and deformations. When an operation loses derivatives, a single Sobolev completion is no longer sufficient, and Hilbert scales and ILH structures enter the picture.

Many identifications that are automatic in finite dimensions cease to be so here. The kinematic tangent space must be distinguished from the operational one, the continuous dual from the algebraic dual, and a strong metric from a weak metric. These differences are more than terminological details: they determine whether a Levi--Civita connection exists, whether the musical map is invertible, and whether a theorem of differential calculus can be applied. The final chapters show how this theory takes concrete form in spaces of maps and connections, gauge groups, diffeomorphisms, and problems such as harmonic maps, Yang--Mills theory, and the Euler equation.

\begin{semblanzaHistorica}{Function spaces as manifolds}
Fréchet opened the way to calculus on abstract spaces, and Banach provided a complete setting in which the fundamental theorems of analysis could be formulated precisely. Later, Eells and Palais developed global tools for infinite-dimensional manifolds. Ebin and Marsden demonstrated the usefulness of Sobolev completions in the study of diffeomorphism groups and fluid equations. With these ideas, many function spaces ceased to be viewed merely as sets of possible solutions and began to be studied as geometric objects in their own right.
\end{semblanzaHistorica}

\chapter{Differential calculus in Banach spaces}\label{cap:calculo-diferencial-banach}

The Sobolev spaces constructed in the preceding parts are themselves spaces in which we can do calculus. A function, a section, or a metric then becomes a point in an infinite-dimensional space, and an energy functional behaves like a function defined on that space. To study its variations, we need a notion of derivative that retains the linear approximation and continuity familiar from the Euclidean case.

The Fréchet derivative provides this language in Banach spaces. It will allow us to formulate chain rules, Taylor expansions, inverse and implicit function theorems, and Lagrange multipliers. In later chapters, these tools will describe manifolds of functions and variational constraints; at the end of the book, critical points of energies defined on Sobolev spaces will again be interpreted as weak solutions of differential equations.

Unless otherwise stated, all normed vector spaces in this chapter are taken over the field of real numbers.

A general reference for the results and exercises on differential
calculus and nonlinear analysis in this chapter is
\cite[Chapters~3, 4, and 6]{DrabekMilota2013}.

\begin{semblanzaHistorica}{Fréchet and Banach: expanding the domain of calculus}
In the early twentieth century, Maurice Fréchet introduced structures capable of measuring convergence beyond Euclidean spaces. Stefan Banach later demonstrated the power of working in complete normed spaces, where approximation and fixed point arguments have limits within the space itself. Fréchet calculus combines these two ideas: the derivative is a continuous linear operator, and completeness allows local estimates to be converted into existence theorems~\cite{Frechet1906,Banach1932}.
\end{semblanzaHistorica}

\section{Fréchet differentiability}
\label{sec:diferenciabilidad-frechet-banach}

The Fréchet derivative generalizes the derivative of functions between Euclidean spaces to the setting of normed spaces. The fundamental idea is the same as in finite dimensions: a function is differentiable at a point if, in a neighborhood of that point, it can be approximated by a continuous linear map so that the approximation error divided by the norm of the increment tends to zero.

The following definition retains the idea of linear approximation but requires the error to be uniformly small in all directions measured by the norm. This uniformity distinguishes the Fréchet derivative from merely directional notions and underpins the fundamental theorems of calculus in infinite dimensions.

For a normed space $(X,\|\cdot\|_X)$, we write $d_X(x,y):=\|x-y\|_X$; thus $B_{d_X}(x,r)$ denotes the ball centered at $x$ with radius $r$ for the metric induced by the norm.

\begin{definition}[Fréchet derivative]
\label{def:derivada-frechet}\glsadd{derivada-frechet}
\index{Frechet derivative@Fréchet derivative}
Let $(X,\|\cdot\|_{X})$ and $(Y,\|\cdot\|_{Y})$ be normed vector spaces, let $U\subseteq X$ be open, and let $f\colon U\longrightarrow Y$. We say that $f$ is \textbf{Fréchet differentiable} at $x_{0}\in U$ if there exists a continuous linear operator $L\in\mathcal{L}(X,Y)$ such that $\displaystyle\lim_{h\to0}\displaystyle\frac{\|f(x_{0}+h)-f(x_{0})-L(h)\|_{Y}}{\|h\|_{X}}=0$. The operator $L$ is called the \textbf{Fréchet derivative} of $f$ at $x_{0}$ and is denoted by $Df(x_{0})$.
\end{definition}

Equivalently, $f$ is Fréchet differentiable at $x_{0}$ if there exist $\rho>0$ and a function $r\colon B_{d_X}(0,\rho)\longrightarrow Y$ such that $B_{d_X}(x_{0},\rho)\subseteq U$, $f(x_{0}+h)=f(x_{0})+Df(x_{0})(h)+r(h)$ for every $h\in B_{d_X}(0,\rho)$, and $\displaystyle\lim_{h\to0}\displaystyle\frac{\|r(h)\|_{Y}}{\|h\|_{X}}=0$.

\begin{proposition}\label{prop:unicidad-derivada-frechet}
The Fréchet derivative of a function, when it exists, is unique.
\end{proposition}

\begin{proof}
Suppose $L_{1},L_{2}\in\mathcal{L}(X,Y)$ satisfy the condition in Definition~\ref{def:derivada-frechet}. Since $U$ is open, there exists $\rho>0$ such that $B_{d_X}(x_{0},\rho)\subseteq U$. Let $v\in X\setminus\{0\}$. For every $t\in\mathbb R$ such that $0<|t|<\displaystyle\frac{\rho}{\|v\|_{X}}$, we have
\[
\frac{\|(L_{1}-L_{2})(v)\|_{Y}}{\|v\|_{X}}
=
\frac{\|(L_{1}-L_{2})(tv)\|_{Y}}{\|tv\|_{X}}
\leq
\frac{\|f(x_{0}+tv)-f(x_{0})-L_{1}(tv)\|_{Y}}{\|tv\|_{X}}
+
\frac{\|f(x_{0}+tv)-f(x_{0})-L_{2}(tv)\|_{Y}}{\|tv\|_{X}}.
\]
As $t\to0$, the right-hand side tends to zero, so $(L_{1}-L_{2})(v)=0$. Equality is also immediate for $v=0$, and therefore $L_{1}=L_{2}$.
\end{proof}

Fréchet differentiability implies, in particular, continuity at the point under consideration.

\begin{proposition}\label{prop:frechet-implica-continuidad}
Let $X$ and $Y$ be normed spaces, let $U\subseteq X$ be open, and let $f\colon U\longrightarrow Y$. If $f$ is Fréchet differentiable at $x_{0}\in U$, then $f$ is continuous at $x_{0}$.
\end{proposition}

\begin{proof}
By the preceding characterization, there exist $\rho>0$ and a function $r\colon B_{d_X}(0,\rho)\longrightarrow Y$ such that $B_{d_X}(x_0,\rho)\subseteq U$, $f(x_0+h)-f(x_0)=Df(x_0)(h)+r(h)$ for every $h\in B_{d_X}(0,\rho)$, and
\[
\lim_{h\to0}\frac{\|r(h)\|_Y}{\|h\|_X}=0.
\]
Let $\varepsilon>0$. There exists $\delta\in(0,\rho)$ such that $\|r(h)\|_Y<\varepsilon\|h\|_X$ whenever $0<\|h\|_X<\delta$. For these increments,
\[
\|f(x_0+h)-f(x_0)\|_Y
\leq
\bigl(\|Df(x_0)\|_{\mathcal L(X,Y)}+\varepsilon\bigr)\|h\|_X.
\]
The right-hand side tends to zero as $h\to0$. Thus $f(x_0+h)\to f(x_0)$, and $f$ is continuous at $x_0$.
\end{proof}

\begin{definition}\label{def:funcion-C1-banach}
\index{C1 function@function of class $C^{1}$}
Let $X$ and $Y$ be normed spaces, let $U\subseteq X$ be open, and let $f\colon U\longrightarrow Y$. We say that $f$ is \textbf{Fréchet differentiable on $U$} if it is Fréchet differentiable at every point of $U$. In this case, the derivative map $Df\colon U\longrightarrow\mathcal{L}(X,Y)$ is defined. We say that $f$ is of \textbf{class $C^{1}$} on $U$ if it is Fréchet differentiable on $U$ and $Df$ is continuous with respect to the operator norm on $\mathcal{L}(X,Y)$.
\end{definition}

The directional derivative studies the behavior of a function only along a fixed direction. Existence of all directional derivatives is generally a weaker condition than Fréchet differentiability.

\begin{definition}[Directional and Gâteaux derivatives]
\label{def:derivada-direccional-gateaux}
\index{directional derivative}
\index{Gateaux derivative@Gâteaux derivative}
Let $(X,\|\cdot\|_{X})$ and $(Y,\|\cdot\|_{Y})$ be normed vector spaces, let $U\subseteq X$ be open, and let $f\colon U\longrightarrow Y$. Given $x_{0}\in U$ and $v\in X$, we say that $f$ is \textbf{differentiable at $x_{0}$ in the direction $v$} if the limit $\displaystyle D_{v}f(x_{0}):=\lim_{t\to0}\displaystyle\frac{f(x_{0}+tv)-f(x_{0})}{t}$ exists. This value is called the \textbf{directional derivative} of $f$ at $x_{0}$ in the direction $v$.

If $D_{v}f(x_{0})$ exists for every $v\in X$ and the map $\delta f(x_{0})\colon X\longrightarrow Y$, defined by $\delta f(x_{0})(v):=D_{v}f(x_{0})$, is linear and continuous, we say that $f$ is \textbf{Gâteaux differentiable} at $x_{0}$. In this case, $\delta f(x_{0})$ is called the \textbf{Gâteaux derivative} of $f$ at $x_{0}$.
\end{definition}

Mere existence of all directional derivatives does not guarantee Gâteaux differentiability, since dependence on the direction need be neither linear nor continuous.

\begin{exercise}\label{ej:derivadas-direccionales-no-gateaux}
Consider the function $f\colon \mathbb{R}^{2}\longrightarrow\mathbb{R}$ defined by
\[
f(x,y)=
\begin{cases}
\displaystyle\frac{x^{2}y}{x^{4}+y^{2}}, & \text{if }(x,y)\neq(0,0),\\[6pt]
0, & \text{if }(x,y)=(0,0).
\end{cases}
\]
Prove that $f$ has a directional derivative at $(0,0)$ in every direction $v=(a,b)\in\mathbb{R}^{2}$ and that
\[
D_{v}f(0,0)=
\begin{cases}
\displaystyle\frac{a^{2}}{b}, & \text{if }b\neq0,\\[6pt]
0, & \text{if }b=0.
\end{cases}
\]
Show that the map $v\longmapsto D_{v}f(0,0)$ is neither linear nor continuous, and conclude that $f$ is not Gâteaux differentiable at $(0,0)$.
\end{exercise}

Fréchet differentiability does guarantee Gâteaux differentiability, and the two derivatives agree.

\begin{proposition}\label{ej:frechet-implica-gateaux}
Let $X$ and $Y$ be normed vector spaces, let $U\subseteq X$ be open, and let $f\colon U\longrightarrow Y$ be Fréchet differentiable at $x_{0}\in U$. Then $f$ is Gâteaux differentiable at $x_{0}$ and
\[
\delta f(x_{0})=Df(x_{0}).
\]
\end{proposition}

\begin{proof}
Set $L:=Df(x_{0})$. By Fréchet differentiability, there exists a function $r$, defined for sufficiently small increments, such that
\[
f(x_{0}+h)-f(x_{0})=Lh+r(h),
\qquad
\frac{\|r(h)\|_{Y}}{\|h\|_{X}}\longrightarrow0
\quad(h\longrightarrow0).
\]
Fix $v\in X$. For sufficiently small $t\neq0$, we have
\[
\frac{f(x_{0}+tv)-f(x_{0})}{t}-Lv
=
\frac{r(tv)}{t}.
\]
If $v=0$, the left-hand side is identically zero. If $v\neq0$, then
\[
\left\|\frac{r(tv)}{t}\right\|_{Y}
=
\|v\|_{X}\frac{\|r(tv)\|_{Y}}{\|tv\|_{X}}
\longrightarrow0
\qquad(t\longrightarrow0).
\]
Consequently, $D_{v}f(x_{0})=Lv$ for every $v\in X$. Since $L$ is linear and continuous, $f$ is Gâteaux differentiable at $x_{0}$ and $\delta f(x_{0})=L=Df(x_{0})$.
\end{proof}

The converse is false, even for scalar functions defined on finite-dimensional Euclidean spaces.

\begin{exercise}\label{ej:gateaux-no-frechet}
Consider the function $g\colon \mathbb{R}^{2}\longrightarrow\mathbb{R}$ defined by
\[
g(x,y)=
\begin{cases}
\displaystyle\frac{x^{3}y}{x^{6}+y^{2}}, & \text{if }(x,y)\neq(0,0),\\[6pt]
0, & \text{if }(x,y)=(0,0).
\end{cases}
\]
Prove that $D_{v}g(0,0)=0$ for every direction $v\in\mathbb{R}^{2}$. Conclude that $g$ is Gâteaux differentiable at $(0,0)$ and that $\delta g(0,0)=0$. Show, however, that $g(t,t^{3})=\displaystyle\frac{1}{2}$ for every $t\neq0$, so $g$ is not continuous at $(0,0)$. Deduce from Proposition~\ref{prop:frechet-implica-continuidad} that $g$ is not Fréchet differentiable at that point.
\end{exercise}

In Hilbert spaces, continuous linear functionals can be identified with vectors by Theorem~\ref{teo:representacion-riesz-hilbert}. This allows the derivative of a scalar functional to be represented by its gradient.

\begin{definition}[Gradient]
\label{def:gradiente-hilbert}
\index{gradient!in a Hilbert space}
Let $(H,\langle\cdot,\cdot\rangle_{H})$ be a real Hilbert space, let $U\subseteq H$ be open, and let $f\colon U\longrightarrow\mathbb{R}$ be Gâteaux differentiable at $x\in U$. By Theorem~\ref{teo:representacion-riesz-hilbert}, there exists a unique vector $\nabla f(x)\in H$ such that $\delta f(x)(v)=\langle v,\nabla f(x)\rangle_{H}$ for every $v\in H$. This vector is called the \textbf{gradient} of $f$ at $x$. If $f$ is Fréchet differentiable at $x$, then $Df(x)(v)=\langle v,\nabla f(x)\rangle_{H}$ for every $v\in H$.
\end{definition}

Before studying mean value formulas, recall that, by Corollary~\ref{cor:caracterizacion-dual-norma}, for every real normed space $Y$ and each $y\in Y\setminus\{0\}$, there exists $\varphi\in Y'$ such that $\|\varphi\|_{Y'}=1$ and $\varphi(y)=\|y\|_{Y}$. In particular, the continuous dual of any normed space separates its points.

The following result estimates the increment of a function in terms of its directional derivatives along the segment joining the points under consideration. The second inequality also gives a quantitative estimate for the error in approximating that increment by the directional derivative at the initial point.

\begin{theorem}[Mean value inequalities]
\label{teo:desigualdades-valor-medio-direccionales}
\index{mean value inequality!in normed spaces}
Let $(X,\|\cdot\|_{X})$ and $(Y,\|\cdot\|_{Y})$ be real normed spaces, let $U\subseteq X$ be open, and let $f\colon U\longrightarrow Y$. Let $x,y\in U$ be such that the segment $[x,y]:=\{x+t(y-x)\mid t\in[0,1]\}$ is contained in $U$. Suppose that, for each $t\in[0,1]$, the directional derivative $D_{y-x}f(x+t(y-x))$ exists. Then
\[
\|f(y)-f(x)\|_{Y}
\leq
\sup_{t\in[0,1]}
\|D_{y-x}f(x+t(y-x))\|_{Y}.
\]
Moreover,
\[
\|f(y)-f(x)-D_{y-x}f(x)\|_{Y}
\leq
\sup_{t\in[0,1]}
\|D_{y-x}f(x+t(y-x))-D_{y-x}f(x)\|_{Y}.
\]

If, in addition, $f$ is Fréchet differentiable at every point of $[x,y]$, then
\[
\|f(y)-f(x)\|_{Y}
\leq
\sup_{t\in[0,1]}
\bigl\|Df(x+t(y-x))(y-x)\bigr\|_{Y}
\leq
\left(
\sup_{t\in[0,1]}
\|Df(x+t(y-x))\|_{\mathcal{L}(X,Y)}
\right)
\|y-x\|_{X},
\]
and
\[
\|f(y)-f(x)-Df(x)(y-x)\|_{Y}
\leq
\sup_{t\in[0,1]}
\bigl\|
\bigl(Df(x+t(y-x))-Df(x)\bigr)(y-x)
\bigr\|_{Y}
\]
\[
\leq
\left(
\sup_{t\in[0,1]}
\|Df(x+t(y-x))-Df(x)\|_{\mathcal{L}(X,Y)}
\right)
\|y-x\|_{X}.
\]

In particular, if $U$ is convex, $f$ is Fréchet differentiable at every point of $U$, and there exists a constant $C\geq0$ such that $\|Df(z)\|_{\mathcal{L}(X,Y)}\leq C$ for every $z\in U$, then $f$ is Lipschitz on $U$ and
\[
\|f(y)-f(x)\|_{Y}
\leq
C\|y-x\|_{X}
\]
for any $x,y\in U$.
\end{theorem}

\begin{proof}
If $f(y)=f(x)$, the first inequality is immediate. Suppose $f(y)\neq f(x)$. By Corollary~\ref{cor:caracterizacion-dual-norma}, there exists $\varphi\in Y'$ such that $\|\varphi\|_{Y'}=1$ and $\varphi(f(y)-f(x))=\|f(y)-f(x)\|_{Y}$.

Define $g\colon [0,1]\longrightarrow\mathbb{R}$ by $g(t):=\varphi(f(x+t(y-x)))$. Existence of $D_{y-x}f(x+t(y-x))$ for each $t\in[0,1]$ implies that $g$ is continuous on $[0,1]$ and differentiable on $(0,1)$. For each $t\in(0,1)$, continuity of $\varphi$ allows us to write
\[
g'(t)
=
\lim_{s\to0}
\varphi\left(
\frac{f(x+(t+s)(y-x))-f(x+t(y-x))}{s}
\right)
=
\varphi\bigl(D_{y-x}f(x+t(y-x))\bigr).
\]
By the classical mean value theorem, there exists $\theta\in(0,1)$ such that $g(1)-g(0)=g'(\theta)$. Consequently,
\[
\|f(y)-f(x)\|_{Y}
=
\varphi(f(y)-f(x))
=
\varphi\bigl(D_{y-x}f(x+\theta(y-x))\bigr)
\]
\[
\leq
\left|
\varphi\bigl(D_{y-x}f(x+\theta(y-x))\bigr)
\right|
\leq
\|D_{y-x}f(x+\theta(y-x))\|_{Y}.
\]
It follows that
\[
\|f(y)-f(x)\|_{Y}
\leq
\sup_{t\in[0,1]}
\|D_{y-x}f(x+t(y-x))\|_{Y}.
\]

To prove the second inequality, the result is immediate if $f(y)-f(x)-D_{y-x}f(x)=0$. Suppose $f(y)-f(x)-D_{y-x}f(x)\neq0$. By Corollary~\ref{cor:caracterizacion-dual-norma}, there exists $\psi\in Y'$ such that $\|\psi\|_{Y'}=1$ and
\[
\psi\bigl(f(y)-f(x)-D_{y-x}f(x)\bigr)
=
\|f(y)-f(x)-D_{y-x}f(x)\|_{Y}.
\]
Define $q\colon [0,1]\longrightarrow\mathbb{R}$ by $q(t):=\psi(f(x+t(y-x))-f(x)-tD_{y-x}f(x))$. The function $q$ is continuous on $[0,1]$ and differentiable on $(0,1)$, and, for each $t\in(0,1)$,
\[
q'(t)
=
\psi\bigl(D_{y-x}f(x+t(y-x))-D_{y-x}f(x)\bigr).
\]
By the classical mean value theorem, there exists $\vartheta\in(0,1)$ such that $q(1)-q(0)=q'(\vartheta)$. Since $q(0)=0$, we obtain
\[
\|f(y)-f(x)-D_{y-x}f(x)\|_{Y}
=
\psi\bigl(D_{y-x}f(x+\vartheta(y-x))-D_{y-x}f(x)\bigr)
\]
\[
\leq
\left|
\psi\bigl(D_{y-x}f(x+\vartheta(y-x))-D_{y-x}f(x)\bigr)
\right|
\leq
\|D_{y-x}f(x+\vartheta(y-x))-D_{y-x}f(x)\|_{Y}.
\]
Therefore,
\[
\|f(y)-f(x)-D_{y-x}f(x)\|_{Y}
\leq
\sup_{t\in[0,1]}
\|D_{y-x}f(x+t(y-x))-D_{y-x}f(x)\|_{Y}.
\]

Now suppose that $f$ is Fréchet differentiable at every point of $[x,y]$. By Proposition~\ref{ej:frechet-implica-gateaux}, for each $t\in[0,1]$ we have $D_{y-x}f(x+t(y-x))=Df(x+t(y-x))(y-x)$. The first inequality already proved then implies that
\[
\|f(y)-f(x)\|_{Y}
\leq
\sup_{t\in[0,1]}
\bigl\|Df(x+t(y-x))(y-x)\bigr\|_{Y}.
\]
On the other hand, for each $t\in[0,1]$,
\[
\bigl\|Df(x+t(y-x))(y-x)\bigr\|_{Y}
\leq
\|Df(x+t(y-x))\|_{\mathcal{L}(X,Y)}
\|y-x\|_{X}.
\]
Taking the supremum over $t\in[0,1]$ gives
\[
\|f(y)-f(x)\|_{Y}
\leq
\left(
\sup_{t\in[0,1]}
\|Df(x+t(y-x))\|_{\mathcal{L}(X,Y)}
\right)
\|y-x\|_{X}.
\]

Likewise, $D_{y-x}f(x)=Df(x)(y-x)$ and, for each $t\in[0,1]$, we have
\[
D_{y-x}f(x+t(y-x))-D_{y-x}f(x)
=
\bigl(Df(x+t(y-x))-Df(x)\bigr)(y-x).
\]
The second inequality already proved implies that
\[
\|f(y)-f(x)-Df(x)(y-x)\|_{Y}
\leq
\sup_{t\in[0,1]}
\bigl\|
\bigl(Df(x+t(y-x))-Df(x)\bigr)(y-x)
\bigr\|_{Y}.
\]
For each $t\in[0,1]$,
\[
\bigl\|
\bigl(Df(x+t(y-x))-Df(x)\bigr)(y-x)
\bigr\|_{Y}
\leq
\|Df(x+t(y-x))-Df(x)\|_{\mathcal{L}(X,Y)}
\|y-x\|_{X}.
\]
Taking the supremum over $t\in[0,1]$ again gives
\[
\|f(y)-f(x)-Df(x)(y-x)\|_{Y}
\leq
\left(
\sup_{t\in[0,1]}
\|Df(x+t(y-x))-Df(x)\|_{\mathcal{L}(X,Y)}
\right)
\|y-x\|_{X}.
\]

Finally, suppose that $U$ is convex and $\|Df(z)\|_{\mathcal{L}(X,Y)}\leq C$ for every $z\in U$. For any $x,y\in U$, the segment $[x,y]$ is contained in $U$ and $\displaystyle\sup_{t\in[0,1]}\|Df(x+t(y-x))\|_{\mathcal{L}(X,Y)}\leq C$. By the first estimate, $\|f(y)-f(x)\|_{Y}\leq C\|y-x\|_{X}$, so $f$ is Lipschitz on $U$.
\end{proof}

The preceding inequalities estimate the increment of the function without requiring continuity of the directional derivatives. Under the additional assumptions that these derivatives vary continuously along the segment and the target space is complete, one can obtain an integral identity for the increment.

\begin{theorem}[Integral mean value formula]
\label{teo: del valor medio}
\label{teo:formula-integral-valor-medio}
\index{integral mean value formula!in Banach spaces}
Let $(X,\|\cdot\|_{X})$ be a real normed space, let $(Y,\|\cdot\|_{Y})$ be a real Banach space, let $U\subseteq X$ be open, and let $f\colon U\longrightarrow Y$. Let $x,y\in U$ be such that the segment $[x,y]:=\{x+t(y-x)\mid t\in[0,1]\}$ is contained in $U$. Suppose $f$ has a directional derivative in the direction $y-x$ at every point of $[x,y]$ and the map $F\colon [0,1]\longrightarrow Y$, defined by $F(t):=D_{y-x}f(x+t(y-x))$, is continuous. Then
\[
f(y)-f(x)
=
\int_{0}^{1}D_{y-x}f(x+t(y-x))\,dt.
\]
\end{theorem}

\begin{proof}
By continuity of $F$ and Theorem~\ref{teo:b5-graves-integral-riemann-banach}, the Riemann integral $\displaystyle\int_{0}^{1}F(t)\,dt$ is well defined. Let $\varphi\in Y'$ and define $g_{\varphi}\colon [0,1]\longrightarrow\mathbb{R}$ by $g_{\varphi}(t):=\varphi(f(x+t(y-x)))$. For each $t\in[0,1]$, interpreting derivatives at the endpoints as one-sided derivatives, continuity of $\varphi$ allows us to write
\[
g_{\varphi}'(t)
=
\lim_{s\to0}
\varphi\left(
\frac{f(x+(t+s)(y-x))-f(x+t(y-x))}{s}
\right)
=
\varphi(F(t)).
\]
Since $F$ is continuous, so is $g_{\varphi}'$. By Theorem~\ref{teo:b5-fundamental-calculo-riemann-banach} and Proposition~\ref{prop:b5-integral-riemann-operadores-lineales},
\[
\varphi(f(y)-f(x))
=
g_{\varphi}(1)-g_{\varphi}(0)
=
\int_{0}^{1}\varphi(F(t))\,dt
=
\varphi\left(\int_{0}^{1}F(t)\,dt\right).
\]
Thus $\displaystyle\varphi\left(f(y)-f(x)-\int_{0}^{1}F(t)\,dt\right)=0$ for every $\varphi\in Y'$. Corollary~\ref{cor:caracterizacion-dual-norma} implies that $\displaystyle f(y)-f(x)-\int_{0}^{1}F(t)\,dt=0$, proving the formula.
\end{proof}

The usual properties of the derivative also hold in this setting.

\begin{theorem}[Properties of the directional derivative]
\label{propiedades derivada direccional}
Let $X$, $Y$, and $Z$ be real normed spaces, let $U\subseteq X$ be an open set, and let $x_0\in U$, $v\in X$, and $f,g\colon U\longrightarrow Y$. If $D_vf(x_0)$ and $D_vg(x_0)$ exist, then:
\begin{enumerate}
 \item $D_v(f+g)(x_0)$ exists and $D_v(f+g)(x_0)=D_vf(x_0)+D_vg(x_0)$.

 \item If $\lambda\in\mathbb R$, then $D_v(\lambda f)(x_0)$ exists and $D_v(\lambda f)(x_0)=\lambda D_vf(x_0)$.

 \item If $B\colon Y\times Y\longrightarrow Z$ is a continuous bilinear map, the function $B(f,g)\colon U\longrightarrow Z$, defined by $B(f,g)(x)=B(f(x),g(x))$, has a directional derivative at $x_0$ in the direction $v$, and
 \[
 D_v\bigl(B(f,g)\bigr)(x_0)
 =
 B\bigl(D_vf(x_0),g(x_0)\bigr)
 +
 B\bigl(f(x_0),D_vg(x_0)\bigr).
 \]

 \item If $Y=\mathbb R$, $g$ is continuous at $x_0$, and $g(x_0)\neq0$, then $D_v\left(\displaystyle\frac{f}{g}\right)(x_0)$ exists and
 \[
 D_v\left(\frac{f}{g}\right)(x_0)
 =
 \frac{D_vf(x_0)g(x_0)-f(x_0)D_vg(x_0)}{g(x_0)^2}.
 \]
\end{enumerate}
\end{theorem}

\begin{proof}
For the first two items, simply apply linearity of limits in the normed space $Y$ term by term.

Let us prove the third. For brevity, write
\[
f_{t}:=f(x_{0}+tv),\qquad g_{t}:=g(x_{0}+tv),
\qquad f_{0}:=f(x_{0}),\qquad g_{0}:=g(x_{0}).
\]
Existence of the directional derivatives implies $f_{t}\to f_{0}$ and $g_{t}\to g_{0}$ as $t\to0$. By bilinearity of $B$,
\[
\frac{B(f_{t},g_{t})-B(f_{0},g_{0})}{t}
=
B\left(\frac{f_{t}-f_{0}}{t},g_{t}\right)
+
B\left(f_{0},\frac{g_{t}-g_{0}}{t}\right).
\]
Continuity of $B$ allows passage to the limit and gives exactly the stated formula.

Finally, suppose $Y=\mathbb R$ and $g(x_{0})\neq0$. Continuity of $g$ at $x_{0}$ ensures that $g_{t}\neq0$ for sufficiently small $t$. Then
\[
\frac{f_{t}/g_{t}-f_{0}/g_{0}}{t}
=
\frac{g_{0}(f_{t}-f_{0})-f_{0}(g_{t}-g_{0})}
     {t\,g_{t}g_{0}}.
\]
Letting $t\to0$ and using $g_{t}\to g_{0}$ and existence of both directional derivatives gives
\[
D_{v}\left(\frac{f}{g}\right)(x_{0})
=
\frac{D_{v}f(x_{0})g(x_{0})-f(x_{0})D_{v}g(x_{0})}
     {g(x_{0})^{2}}.
\]
\end{proof}

\begin{theorem}[Chain rule for directional derivatives]
\label{regla de la cadena derivadas direccionales}
Let $X$, $Y$, and $Z$ be real normed spaces, let $U\subseteq X$ and $V\subseteq Y$ be open sets, and let $f\colon U\longrightarrow Y$ and $g\colon V\longrightarrow Z$ be functions such that $f(U)\subseteq V$. Let $x_{0}\in U$ and $v\in X$. If $g$ is Fréchet differentiable at $f(x_{0})$ and $D_{v}f(x_{0})$ exists, then $D_{v}(g\circ f)(x_{0})$ exists and
\[
D_{v}(g\circ f)(x_{0})
=
Dg(f(x_{0}))\bigl(D_{v}f(x_{0})\bigr).
\]
\end{theorem}

\begin{proof}
Since $U$ and $V$ are open, there exist $\delta_{0}>0$ and $\rho>0$ such that $x_{0}+tv\in U$ for every $t\in\mathbb R$ with $|t|<\delta_{0}$ and $B_{d_Y}(f(x_{0}),\rho)\subseteq V$. We proceed by cases.

\begin{enumerate}
 \item Case $1$: $D_{v}f(x_{0})\neq0$.

 Let $\varepsilon>0$. Since $g$ is Fréchet differentiable at $f(x_{0})$, there exists $\delta_{1}\in(0,\rho)$ such that if $0<\|h\|_{Y}<\delta_{1}$, then
 \[
 \frac{\|g(f(x_{0})+h)-g(f(x_{0}))-Dg(f(x_{0}))(h)\|_{Z}}{\|h\|_{Y}}
 <
 \frac{\varepsilon}{2\bigl(1+\|D_{v}f(x_{0})\|_{Y}\bigr)}.
 \]

 On the other hand, since $D_{v}f(x_{0})$ exists, there exists $\delta_{2}>0$ such that if $0<|t|<\delta_{2}$, then
 \[
 \left\|
 \frac{f(x_{0}+tv)-f(x_{0})}{t}-D_{v}f(x_{0})
 \right\|_{Y}
 <
 \min\left\{
 1,
 \frac{\|D_{v}f(x_{0})\|_{Y}}{2},
 \frac{\varepsilon}{2\bigl(1+\|Dg(f(x_{0}))\|_{\mathcal L(Y,Z)}\bigr)}
 \right\}.
 \]

 Let $\displaystyle\delta=\min\left\{\delta_{0},\delta_{2},\displaystyle\frac{\delta_{1}}{1+\|D_{v}f(x_{0})\|_{Y}}\right\}>0$ and let $t\in\mathbb R$ be such that $0<|t|<\delta$. First note that
 \[
 \frac{\|f(x_{0}+tv)-f(x_{0})\|_{Y}}{|t|}
 \leq
 \left\|
 \frac{f(x_{0}+tv)-f(x_{0})}{t}-D_{v}f(x_{0})
 \right\|_{Y}
 +
 \|D_{v}f(x_{0})\|_{Y}
 <
 1+\|D_{v}f(x_{0})\|_{Y}.
 \]

 On the other hand, by the reverse triangle inequality,
 \[
 \frac{\|f(x_{0}+tv)-f(x_{0})\|_{Y}}{|t|}
 \geq
 \|D_{v}f(x_{0})\|_{Y}
 -
 \left\|
 \frac{f(x_{0}+tv)-f(x_{0})}{t}-D_{v}f(x_{0})
 \right\|_{Y}
 >
 \frac{\|D_{v}f(x_{0})\|_{Y}}{2}>0.
 \]

 Consequently, $f(x_{0}+tv)\neq f(x_{0})$ and $0<\|f(x_{0}+tv)-f(x_{0})\|_{Y}<|t|\bigl(1+\|D_{v}f(x_{0})\|_{Y}\bigr)<\delta_{1}$.

 Adding and subtracting $\displaystyle\frac{Dg(f(x_{0}))\bigl(f(x_{0}+tv)-f(x_{0})\bigr)}{t}$ and using the triangle inequality, the quantity
 \[
 \left\|
 \frac{g(f(x_{0}+tv))-g(f(x_{0}))}{t}
 -
 Dg(f(x_{0}))\bigl(D_{v}f(x_{0})\bigr)
 \right\|_{Z}
 \]
 is bounded above by the sum of
 \[
 \frac{\|g(f(x_{0}+tv))-g(f(x_{0}))-Dg(f(x_{0}))(f(x_{0}+tv)-f(x_{0}))\|_{Z}}{|t|}
 \]
 and
 \[
 \|Dg(f(x_{0}))\|_{\mathcal L(Y,Z)}
 \left\|
 \frac{f(x_{0}+tv)-f(x_{0})}{t}-D_{v}f(x_{0})
 \right\|_{Y}.
 \]

 The first summand equals
 \[
 \frac{\|g(f(x_{0}+tv))-g(f(x_{0}))-Dg(f(x_{0}))(f(x_{0}+tv)-f(x_{0}))\|_{Z}}
 {\|f(x_{0}+tv)-f(x_{0})\|_{Y}}
 \frac{\|f(x_{0}+tv)-f(x_{0})\|_{Y}}{|t|},
 \]
 so
 \[
 \frac{\|g(f(x_{0}+tv))-g(f(x_{0}))-Dg(f(x_{0}))(f(x_{0}+tv)-f(x_{0}))\|_{Z}}{|t|}
 <
 \frac{\varepsilon}{2}.
 \]

 Likewise,
 \[
 \|Dg(f(x_{0}))\|_{\mathcal L(Y,Z)}
 \left\|
 \frac{f(x_{0}+tv)-f(x_{0})}{t}-D_{v}f(x_{0})
 \right\|_{Y}
 \leq
 \frac{\varepsilon\|Dg(f(x_{0}))\|_{\mathcal L(Y,Z)}}
 {2\bigl(1+\|Dg(f(x_{0}))\|_{\mathcal L(Y,Z)}\bigr)}
 <
 \frac{\varepsilon}{2}.
 \]

 Therefore,
 \[
 \left\|
 \frac{g(f(x_{0}+tv))-g(f(x_{0}))}{t}
 -
 Dg(f(x_{0}))\bigl(D_{v}f(x_{0})\bigr)
 \right\|_{Z}
 <
 \varepsilon.
 \]
 Consequently, $D_{v}(g\circ f)(x_{0})=Dg(f(x_{0}))\bigl(D_{v}f(x_{0})\bigr)$.

 \item Case $2$: $D_{v}f(x_{0})=0$.

 In this case, we must prove that $D_{v}(g\circ f)(x_{0})=0$. Let $\varepsilon>0$. Since $g$ is Fréchet differentiable at $f(x_{0})$, there exists $\delta_{1}\in(0,\rho)$ such that if $0<\|h\|_{Y}<\delta_{1}$, then
 \[
 \frac{\|g(f(x_{0})+h)-g(f(x_{0}))-Dg(f(x_{0}))(h)\|_{Z}}{\|h\|_{Y}}<1.
 \]

 On the other hand, since $D_{v}f(x_{0})=0$, there exists $\delta_{2}>0$ such that if $0<|t|<\delta_{2}$, then
 \[
 \frac{\|f(x_{0}+tv)-f(x_{0})\|_{Y}}{|t|}
 <
 \min\left\{
 1,
 \frac{\varepsilon}{1+\|Dg(f(x_{0}))\|_{\mathcal L(Y,Z)}}
 \right\}.
 \]

 Let $\displaystyle \delta=\min\{\delta_{0},\delta_{1},\delta_{2}\}>0$ and let $t\in\mathbb R$ be such that $0<|t|<\delta$. If $f(x_{0}+tv)=f(x_{0})$, then
 \[
 \left\|
 \frac{g(f(x_{0}+tv))-g(f(x_{0}))}{t}
 \right\|_{Z}
 =
 0
 <
 \varepsilon.
 \]

 Now suppose $f(x_{0}+tv)\neq f(x_{0})$. In this case, $0<\|f(x_{0}+tv)-f(x_{0})\|_{Y}<|t|<\delta\leq\delta_{1}$.

 Adding and subtracting $\displaystyle\frac{Dg(f(x_{0}))\bigl(f(x_{0}+tv)-f(x_{0})\bigr)}{t}$ and using the triangle inequality, the quantity
 \[
 \left\|
 \frac{g(f(x_{0}+tv))-g(f(x_{0}))}{t}
 \right\|_{Z}
 \]
 is bounded above by the sum of
 \[
 \frac{\|g(f(x_{0}+tv))-g(f(x_{0}))-Dg(f(x_{0}))(f(x_{0}+tv)-f(x_{0}))\|_{Z}}{|t|}
 \]
 and
 \[
 \|Dg(f(x_{0}))\|_{\mathcal L(Y,Z)}
 \frac{\|f(x_{0}+tv)-f(x_{0})\|_{Y}}{|t|}.
 \]

 The first summand equals
 \[
 \frac{\|g(f(x_{0}+tv))-g(f(x_{0}))-Dg(f(x_{0}))(f(x_{0}+tv)-f(x_{0}))\|_{Z}}
 {\|f(x_{0}+tv)-f(x_{0})\|_{Y}}
 \frac{\|f(x_{0}+tv)-f(x_{0})\|_{Y}}{|t|},
 \]
 so
 \[
 \frac{\|g(f(x_{0}+tv))-g(f(x_{0}))-Dg(f(x_{0}))(f(x_{0}+tv)-f(x_{0}))\|_{Z}}{|t|}
 <
 \frac{\varepsilon}{1+\|Dg(f(x_{0}))\|_{\mathcal L(Y,Z)}}.
 \]

 Likewise,
 \[
 \|Dg(f(x_{0}))\|_{\mathcal L(Y,Z)}
 \frac{\|f(x_{0}+tv)-f(x_{0})\|_{Y}}{|t|}
 \leq
 \frac{\varepsilon\|Dg(f(x_{0}))\|_{\mathcal L(Y,Z)}}
 {1+\|Dg(f(x_{0}))\|_{\mathcal L(Y,Z)}}.
 \]

 Therefore,
 \[
 \left\|
 \frac{g(f(x_{0}+tv))-g(f(x_{0}))}{t}
 \right\|_{Z}
 <
 \varepsilon.
 \]
 Consequently, $D_{v}(g\circ f)(x_{0})=0=Dg(f(x_{0}))\bigl(D_{v}f(x_{0})\bigr)$.
\end{enumerate}
\end{proof}

The chain rule for directional derivatives immediately gives the corresponding rule for the Gâteaux derivative. The version for the Fréchet derivative also requires checking that the linear approximation is uniform in the increments.

\begin{corollary}[Chain rule for Gâteaux and Fréchet derivatives]
\label{cor:regla-cadena-gateaux-frechet}
\index{chain rule!for the Gâteaux derivative}
\index{chain rule!for the Fréchet derivative}
Let $X$, $Y$, and $Z$ be normed spaces, let $U\subseteq X$ and $V\subseteq Y$ be open sets, and let $f\colon U\longrightarrow Y$ and $g\colon V\longrightarrow Z$ be functions such that $f(U)\subseteq V$. Let $x_{0}\in U$ and suppose $g$ is Fréchet differentiable at $f(x_{0})$.

If $f$ is Gâteaux differentiable at $x_{0}$, then $g\circ f$ is Gâteaux differentiable at $x_{0}$ and
\[
\delta(g\circ f)(x_{0})
=
Dg(f(x_{0}))\circ\delta f(x_{0}).
\]

If $f$ is Fréchet differentiable at $x_{0}$, then $g\circ f$ is Fréchet differentiable at $x_{0}$ and
\[
D(g\circ f)(x_{0})
=
Dg(f(x_{0}))\circ Df(x_{0}).
\]
\end{corollary}

\begin{proof}
First suppose $f$ is Gâteaux differentiable at $x_{0}$. For each $v\in X$, Theorem~\ref{regla de la cadena derivadas direccionales} implies that
\[
D_{v}(g\circ f)(x_{0})
=
Dg(f(x_{0}))\bigl(D_{v}f(x_{0})\bigr)
=
Dg(f(x_{0}))\bigl(\delta f(x_{0})(v)\bigr).
\]
Since $Dg(f(x_{0}))\circ\delta f(x_{0})\in\mathcal L(X,Z)$, we conclude that $g\circ f$ is Gâteaux differentiable at $x_{0}$ and $\delta(g\circ f)(x_{0})=Dg(f(x_{0}))\circ\delta f(x_{0})$.

Now suppose $f$ is Fréchet differentiable at $x_{0}$. By Proposition~\ref{ej:frechet-implica-gateaux}, the first part shows that $g\circ f$ is Gâteaux differentiable at $x_{0}$ and the only candidate for its Fréchet derivative is $Dg(f(x_{0}))\circ Df(x_{0})$.

For each $h\in X\setminus\{0\}$ such that $x_{0}+h\in U$, adding and subtracting $Dg(f(x_{0}))(f(x_{0}+h)-f(x_{0}))$ and using the triangle inequality, the quantity
\[
\frac{
\bigl\|
g(f(x_{0}+h))-g(f(x_{0}))
-Dg(f(x_{0}))\bigl(Df(x_{0})(h)\bigr)
\bigr\|_{Z}
}{
\|h\|_{X}
}
\]
is bounded above by the sum of
\[
\frac{
\bigl\|
g(f(x_{0}+h))-g(f(x_{0}))
-Dg(f(x_{0}))\bigl(f(x_{0}+h)-f(x_{0})\bigr)
\bigr\|_{Z}
}{
\|h\|_{X}
}
\]
and
\[
\|Dg(f(x_{0}))\|_{\mathcal L(Y,Z)}
\frac{
\|f(x_{0}+h)-f(x_{0})-Df(x_{0})(h)\|_{Y}
}{
\|h\|_{X}
}.
\]

If $f(x_{0}+h)\neq f(x_{0})$, the first quotient can be written as
\[
\frac{
\bigl\|
g(f(x_{0}+h))-g(f(x_{0}))
-Dg(f(x_{0}))\bigl(f(x_{0}+h)-f(x_{0})\bigr)
\bigr\|_{Z}
}{
\|f(x_{0}+h)-f(x_{0})\|_{Y}
}
\frac{
\|f(x_{0}+h)-f(x_{0})\|_{Y}
}{
\|h\|_{X}
}.
\]
If $f(x_{0}+h)=f(x_{0})$, that quotient equals zero. Moreover,
\[
\frac{\|f(x_{0}+h)-f(x_{0})\|_{Y}}{\|h\|_{X}}
\leq
\frac{
\|f(x_{0}+h)-f(x_{0})-Df(x_{0})(h)\|_{Y}
}{
\|h\|_{X}
}
+
\|Df(x_{0})\|_{\mathcal L(X,Y)}.
\]
By Fréchet differentiability of $f$ at $x_{0}$, this last quotient remains bounded as $h\to0$, while $f(x_{0}+h)\to f(x_{0})$ by Proposition~\ref{prop:frechet-implica-continuidad}. Fréchet differentiability of $g$ at $f(x_{0})$ then implies that the first summand tends to zero, and differentiability of $f$ implies the same for the second. Thus
\[
\lim_{h\to0}
\frac{
\bigl\|
g(f(x_{0}+h))-g(f(x_{0}))
-Dg(f(x_{0}))\bigl(Df(x_{0})(h)\bigr)
\bigr\|_{Z}
}{
\|h\|_{X}
}
=
0.
\]
Consequently, $g\circ f$ is Fréchet differentiable at $x_{0}$ and $D(g\circ f)(x_{0})=Dg(f(x_{0}))\circ Df(x_{0})$.
\end{proof}

\begin{corollary}[Composition with a continuous linear operator]
\label{cor:composicion-operador-lineal-continuo}
\index{derivative!composition with a continuous linear operator}
Let $X$, $Y$, and $Z$ be normed spaces, let $U\subseteq X$ be open, let $f\colon U\longrightarrow Y$, and let $A\in\mathcal L(Y,Z)$. Given $x_{0}\in U$ and $v\in X$, the following assertions hold:
\begin{enumerate}
 \item If $D_{v}f(x_{0})$ exists, then $D_{v}(A\circ f)(x_{0})$ exists and $D_{v}(A\circ f)(x_{0})=A\bigl(D_{v}f(x_{0})\bigr)$.

 \item If $f$ is Gâteaux differentiable at $x_{0}$, then $A\circ f$ is Gâteaux differentiable at $x_{0}$ and $\delta(A\circ f)(x_{0})=A\circ\delta f(x_{0})$.

 \item If $f$ is Fréchet differentiable at $x_{0}$, then $A\circ f$ is Fréchet differentiable at $x_{0}$ and $D(A\circ f)(x_{0})=A\circ Df(x_{0})$.
\end{enumerate}
\end{corollary}

\begin{proof}
For each $y\in Y$ and each $h\in Y\setminus\{0\}$, we have
\[
\frac{\|A(y+h)-A(y)-A(h)\|_{Z}}{\|h\|_{Y}}=0.
\]
Thus $A$ is Fréchet differentiable at every point of $Y$ and $DA(y)=A$. The three assertions follow from Theorem~\ref{regla de la cadena derivadas direccionales} and Corollary~\ref{cor:regla-cadena-gateaux-frechet}.
\end{proof}

Existence of the Gâteaux derivative alone does not guarantee Fréchet differentiability. However, continuity of the Gâteaux derivative in the operator norm provides a sufficient condition.

\begin{proposition}[Continuity of the Gâteaux derivative]
\label{prop:continuidad-gateaux-implica-frechet}
\index{Gateaux derivative@Gâteaux derivative!continuity}
\index{Frechet derivative@Fréchet derivative!criterion using the Gâteaux derivative}
Let $X$ and $Y$ be normed spaces, let $U\subseteq X$ be open, let $x_{0}\in U$, and let $f\colon U\longrightarrow Y$. Suppose there exists $r>0$ such that $B_{d_X}(x_{0},r)\subseteq U$, $f$ is Gâteaux differentiable at every point of $B_{d_X}(x_{0},r)$, and the map $\delta f\colon B_{d_X}(x_{0},r)\longrightarrow\mathcal L(X,Y)$ is continuous at $x_{0}$. Then $f$ is Fréchet differentiable at $x_{0}$ and $Df(x_{0})=\delta f(x_{0})$.
\end{proposition}

\begin{proof}
Let $\varepsilon>0$. By continuity of $\delta f$ at $x_{0}$, there exists $\delta\in(0,r)$ such that $\|\delta f(z)-\delta f(x_{0})\|_{\mathcal L(X,Y)}<\varepsilon$ for every $z\in B_{d_X}(x_{0},\delta)$.

Let $h\in X$ be such that $0<\|h\|_{X}<\delta$. For each $t\in[0,1]$, we have $\|x_{0}+th-x_{0}\|_{X}=t\|h\|_{X}<\delta$, so $x_{0}+th\in B_{d_X}(x_{0},\delta)$.

Applying the second estimate in Theorem~\ref{teo:desigualdades-valor-medio-direccionales} to the points $x_{0}$ and $x_{0}+h$ gives
\[
\|f(x_{0}+h)-f(x_{0})-D_{h}f(x_{0})\|_{Y}
\leq
\sup_{t\in[0,1]}
\|D_{h}f(x_{0}+th)-D_{h}f(x_{0})\|_{Y}.
\]
Since $f$ is Gâteaux differentiable at every point of the segment, for each $t\in[0,1]$ we have $D_{h}f(x_{0}+th)=\delta f(x_{0}+th)(h)$ and $D_{h}f(x_{0})=\delta f(x_{0})(h)$. Thus
\[
\|f(x_{0}+h)-f(x_{0})-\delta f(x_{0})(h)\|_{Y}
\leq
\sup_{t\in[0,1]}
\bigl\|
\bigl(\delta f(x_{0}+th)-\delta f(x_{0})\bigr)(h)
\bigr\|_{Y}.
\]
By the definition of the operator norm and the choice of $\delta$, for each $t\in[0,1]$ we have
\[
\bigl\|
\bigl(\delta f(x_{0}+th)-\delta f(x_{0})\bigr)(h)
\bigr\|_{Y}
\leq
\|\delta f(x_{0}+th)-\delta f(x_{0})\|_{\mathcal L(X,Y)}
\|h\|_{X}
<
\varepsilon\|h\|_{X}.
\]
Consequently, $\|f(x_{0}+h)-f(x_{0})-\delta f(x_{0})(h)\|_{Y}\leq\varepsilon\|h\|_{X}$ and therefore
\[
\frac{
\|f(x_{0}+h)-f(x_{0})-\delta f(x_{0})(h)\|_{Y}
}{
\|h\|_{X}
}
\leq
\varepsilon.
\]
This proves that $f$ is Fréchet differentiable at $x_{0}$ and that $Df(x_{0})=\delta f(x_{0})$.
\end{proof}

The preceding proposition allows us to characterize functions of class $C^{1}$ using the Gâteaux derivative and the notion of differentiability already defined.

\begin{corollary}[Characterization of functions of class \(C^{1}\)]
\label{cor:caracterizacion-C1-gateaux}
\index{function of class \(C^{1}\)!characterization using Gâteaux derivatives}
Let $X$ and $Y$ be normed spaces, let $U\subseteq X$ be open, and let $f\colon U\longrightarrow Y$. Then $f\in C^{1}(U,Y)$ if and only if $f$ is Gâteaux differentiable at every point of $U$ and the map $\delta f\colon U\longrightarrow\mathcal L(X,Y)$ is continuous. In this case, $Df(x)=\delta f(x)$ for every $x\in U$.
\end{corollary}

\begin{proof}
If $f\in C^{1}(U,Y)$, Proposition~\ref{ej:frechet-implica-gateaux} implies that $f$ is Gâteaux differentiable at every point of $U$ and that $\delta f=Df$. Continuity of $\delta f$ follows from continuity of $Df$.

Conversely, suppose $f$ is Gâteaux differentiable at every point of $U$ and $\delta f\colon U\longrightarrow\mathcal L(X,Y)$ is continuous. Let $x\in U$. Since $U$ is open, there exists $r>0$ such that $B_{d_X}(x,r)\subseteq U$. Proposition~\ref{prop:continuidad-gateaux-implica-frechet} implies that $f$ is Fréchet differentiable at $x$ and that $Df(x)=\delta f(x)$. Since $x\in U$ was arbitrary, $f$ is Fréchet differentiable on $U$ and $Df=\delta f$. Continuity of $\delta f$ then implies continuity of $Df$. By Definition~\ref{def:funcion-C1-banach}, we conclude that $f\in C^{1}(U,Y)$.
\end{proof}
The notion of partial derivative also extends naturally to normed spaces. Here it is useful to distinguish Gâteaux and Fréchet partial derivatives according to the notion of differentiability used when the other variables are held fixed. On the product space $X\times Y$, we use the norm
$\|(h,k)\|_{X\times Y}:=\|h\|_{X}+\|k\|_{Y}$.

\begin{definition}[Gâteaux and Fréchet partial derivatives]
\label{def:derivadas-parciales-gateaux-frechet}
\index{partial derivative!Gâteaux}
\index{partial derivative!Fréchet}
Let $X$, $Y$, and $Z$ be normed spaces, let $U\subseteq X$ and $V\subseteq Y$ be open, let $f\colon U\times V\longrightarrow Z$, and let $(x,y)\in U\times V$.

If the function $f(\,\cdot\,,y)\colon U\longrightarrow Z$ is Gâteaux differentiable at $x$, the \textbf{Gâteaux partial derivative of $f$ with respect to the first variable} is defined by
$\delta_{1}f(x,y):=\delta(f(\,\cdot\,,y))(x)\in\mathcal L(X,Z)$. If $f(x,\,\cdot\,)\colon V\longrightarrow Z$ is Gâteaux differentiable at $y$, the \textbf{Gâteaux partial derivative with respect to the second variable} is defined by
$\delta_{2}f(x,y):=\delta(f(x,\,\cdot\,))(y)\in\mathcal L(Y,Z)$.

If $f(\,\cdot\,,y)$ is Fréchet differentiable at $x$, the \textbf{Fréchet partial derivative of $f$ with respect to the first variable} is defined by
$D_{1}f(x,y):=D(f(\,\cdot\,,y))(x)\in\mathcal L(X,Z)$. If $f(x,\,\cdot\,)$ is Fréchet differentiable at $y$, the \textbf{Fréchet partial derivative with respect to the second variable} is defined by
$D_{2}f(x,y):=D(f(x,\,\cdot\,))(y)\in\mathcal L(Y,Z)$.
\end{definition}

Partial derivatives are obtained by restricting the total derivative to the directions corresponding to each factor of the product.

\begin{proposition}
\label{prop:derivada-total-implica-parciales}
\index{partial derivative!relation to the total derivative}
Let $X$, $Y$, and $Z$ be normed spaces, let $U\subseteq X$ and $V\subseteq Y$ be open, let $f\colon U\times V\longrightarrow Z$, and let $(x,y)\in U\times V$.

If $f$ is Gâteaux differentiable at $(x,y)$, then $\delta_{1}f(x,y)$ and $\delta_{2}f(x,y)$ exist and, for every $(h,k)\in X\times Y$,
\begin{equation}
\label{eq:descomposicion-derivada-total-gateaux}
\delta f(x,y)(h,k)
=
\delta_{1}f(x,y)(h)
+
\delta_{2}f(x,y)(k).
\end{equation}

If $f$ is Fréchet differentiable at $(x,y)$, then $D_{1}f(x,y)$ and $D_{2}f(x,y)$ exist and, for every $(h,k)\in X\times Y$,
\begin{equation}
\label{eq:descomposicion-derivada-total-frechet}
Df(x,y)(h,k)
=
D_{1}f(x,y)(h)
+
D_{2}f(x,y)(k).
\end{equation}
\end{proposition}

\begin{proof}
First suppose $f$ is Gâteaux differentiable at $(x,y)$. Let $h\in X$. When only the first component varies, the increment $th$ of $x$ corresponds to the increment $(th,0)=t(h,0)$ of the point $(x,y)$. Thus
\[
\delta_{1}f(x,y)(h)
=
\lim_{t\to0}
\frac{f(x+th,y)-f(x,y)}{t}
=
\lim_{t\to0}
\frac{f\bigl((x,y)+t(h,0)\bigr)-f(x,y)}{t}
=
\delta f(x,y)(h,0).
\]
This proves that $\delta_{1}f(x,y)$ exists and agrees with the restriction of $\delta f(x,y)$ to the subspace $X\times\{0\}$. For the second variable, varying only that component corresponds to the increment $(0,tk)=t(0,k)$. Hence, for each $k\in Y$,
\[
\delta_{2}f(x,y)(k)
=
\lim_{t\to0}
\frac{f(x,y+tk)-f(x,y)}{t}
=
\delta f(x,y)(0,k).
\]
Since $(h,k)=(h,0)+(0,k)$ and $\delta f(x,y)$ is linear, it follows that
\[
\delta f(x,y)(h,k)
=
\delta f(x,y)(h,0)+\delta f(x,y)(0,k)
=
\delta_{1}f(x,y)(h)+\delta_{2}f(x,y)(k),
\]
which proves \eqref{eq:descomposicion-derivada-total-gateaux}.

Now suppose $f$ is Fréchet differentiable at $(x,y)$. For each $h\in X\setminus\{0\}$ such that $x+h\in U$, we have
\[
\frac{
\|f(x+h,y)-f(x,y)-Df(x,y)(h,0)\|_{Z}
}{
\|h\|_{X}
}
=
\frac{
\|f\bigl((x,y)+(h,0)\bigr)-f(x,y)-Df(x,y)(h,0)\|_{Z}
}{
\|(h,0)\|_{X\times Y}
}.
\]
The right-hand side tends to zero as $h\to0$, since $(h,0)\to(0,0)$ in $X\times Y$. Consequently, $f(\,\cdot\,,y)$ is Fréchet differentiable at $x$ and
$D_{1}f(x,y)(h)=Df(x,y)(h,0)$ for every $h\in X$.

For the second variable, if $k\in Y\setminus\{0\}$ and $y+k\in V$, then
\[
\frac{\|f(x,y+k)-f(x,y)-Df(x,y)(0,k)\|_{Z}}{\|k\|_{Y}}
=
\frac{\|f\bigl((x,y)+(0,k)\bigr)-f(x,y)-Df(x,y)(0,k)\|_{Z}}
{\|(0,k)\|_{X\times Y}}.
\]
The last quotient tends to zero as $k\to0$. Thus
$D_{2}f(x,y)(k)=Df(x,y)(0,k)$ for every $k\in Y$. Linearity of $Df(x,y)$ then implies
\[
Df(x,y)(h,k)
=
Df(x,y)(h,0)+Df(x,y)(0,k)
=
D_{1}f(x,y)(h)+D_{2}f(x,y)(k),
\]
which proves \eqref{eq:descomposicion-derivada-total-frechet}.
\end{proof}

The converse generally does not follow from mere existence of the partial derivatives. However, continuity of one of them in the operator norm allows us to recover total differentiability.

\begin{proposition}[Differentiability from partial derivatives]
\label{prop:diferenciabilidad-mediante-derivadas-parciales}
\index{partial derivative!differentiability criterion}
Let $X$, $Y$, and $Z$ be normed spaces, let $U\subseteq X$ and $V\subseteq Y$ be open, let $f\colon U\times V\longrightarrow Z$, and let $(x_{0},y_{0})\in U\times V$. Suppose $\delta_{1}f(x_{0},y_{0})$ exists, $\delta_{2}f$ is defined on an open neighborhood $W\subseteq U\times V$ of $(x_{0},y_{0})$, and the map
\[
(x,y)\longmapsto\delta_{2}f(x,y),
\qquad
W\longrightarrow\mathcal L(Y,Z),
\]
is continuous at $(x_{0},y_{0})$. Then $f$ is Gâteaux differentiable at $(x_{0},y_{0})$, and identity \eqref{eq:descomposicion-derivada-total-gateaux} holds with $(x,y)=(x_{0},y_{0})$.

If, in addition, $D_{1}f(x_{0},y_{0})$ exists, then $f$ is Fréchet differentiable at $(x_{0},y_{0})$, and identity \eqref{eq:descomposicion-derivada-total-frechet} holds with $(x,y)=(x_{0},y_{0})$. In this case,
\[
D_{2}f(x_{0},y_{0})
=
\delta_{2}f(x_{0},y_{0}).
\]
\end{proposition}

\begin{proof}
Let $W\subseteq U\times V$ be the open neighborhood in the statement. Take $r>0$ such that $B_{d_{X\times Y}}((x_0,y_0),r)\subseteq W$.

Fix $(h,k)\in X\times Y$. The case $(h,k)=(0,0)$ is immediate. Suppose $(h,k)\neq(0,0)$ and
$\displaystyle 0<|t|<\displaystyle\frac{r}{\|h\|_X+\|k\|_Y}$. For each $\tau\in[0,1]$, we have $\|(th,\tau tk)\|_{X\times Y}\leq |t|(\|h\|_X+\|k\|_Y)<r$; hence the segment joining $(x_0+th,y_0)$ to $(x_0+th,y_0+tk)$ is contained in $W$.

Adding and subtracting $f(x_0+th,y_0)$ gives
\[
\frac{f(x_0+th,y_0+tk)-f(x_0,y_0)}{t}
-\delta_1f(x_0,y_0)(h)-\delta_2f(x_0,y_0)(k)\] \[=
\frac{f(x_0+th,y_0+tk)-f(x_0+th,y_0)-\delta_2f(x_0,y_0)(tk)}{t}+
\frac{f(x_0+th,y_0)-f(x_0,y_0)}{t}-\delta_1f(x_0,y_0)(h).
\]

To estimate the first term, apply the first mean value inequality in Theorem~\ref{teo:desigualdades-valor-medio-direccionales} to the map $v\longmapsto f(x_0+th,v)-\delta_2f(x_0,y_0)(v-y_0)$ between the points $y_0$ and $y_0+tk$. Its directional derivative in the direction $tk$ is $\bigl(\delta_2f(x_0+th,v)-\delta_2f(x_0,y_0)\bigr)(tk)$. Dividing the resulting estimate by $|t|$ gives
\[
\begin{aligned}
&\left\|
\frac{f(x_0+th,y_0+tk)-f(x_0+th,y_0)-\delta_2f(x_0,y_0)(tk)}{t}
\right\|_Z\\
&\quad\leq\displaystyle\sup_{\tau\in[0,1]}
\|\delta_2f(x_0+th,y_0+\tau tk)-\delta_2f(x_0,y_0)\|_{\mathcal L(Y,Z)}\|k\|_Y.
\end{aligned}
\]

Let $\varepsilon>0$. By continuity of $\delta_2f$ at $(x_0,y_0)$, there exists $\eta\in(0,r)$ such that
$\|\delta_2f(u,v)-\delta_2f(x_0,y_0)\|_{\mathcal L(Y,Z)}<\varepsilon$ whenever $(u,v)\in W$ and $\|(u-x_0,v-y_0)\|_{X\times Y}<\eta$. If
$\displaystyle 0<|t|<\displaystyle\frac{\eta}{\|h\|_X+\|k\|_Y}$, then $\|(th,\tau tk)\|_{X\times Y}<\eta$ for every $\tau\in[0,1]$, and consequently
\[
\left\|
\frac{f(x_0+th,y_0+tk)-f(x_0+th,y_0)}{t}
-\delta_2f(x_0,y_0)(k)
\right\|_Z
\leq\varepsilon\|k\|_Y.
\]
Thus
\[
\lim_{t\to0}
\left(
\frac{f(x_0+th,y_0+tk)-f(x_0+th,y_0)}{t}
-\delta_2f(x_0,y_0)(k)
\right)=0.
\]
On the other hand, existence of $\delta_1f(x_0,y_0)$ implies
\[
\lim_{t\to0}
\left(
\frac{f(x_0+th,y_0)-f(x_0,y_0)}{t}
-\delta_1f(x_0,y_0)(h)
\right)=0.
\]
Combining the two limits, we conclude that the directional derivative of $f$ at $(x_0,y_0)$ in the direction $(h,k)$ exists and equals
\[
\delta_1f(x_0,y_0)(h)+\delta_2f(x_0,y_0)(k).
\]
The map $(h,k)\longmapsto\delta_1f(x_0,y_0)(h)+\delta_2f(x_0,y_0)(k)$ is linear. Moreover,
\[
\|\delta_1f(x_0,y_0)(h)+\delta_2f(x_0,y_0)(k)\|_Z\leq
\max\bigl\{\|\delta_1f(x_0,y_0)\|_{\mathcal L(X,Z)},
\|\delta_2f(x_0,y_0)\|_{\mathcal L(Y,Z)}\bigr\}
\|(h,k)\|_{X\times Y},
\]
so it is continuous. Thus $f$ is Gâteaux differentiable at $(x_0,y_0)$ and \eqref{eq:descomposicion-derivada-total-gateaux} holds.

Now suppose $D_1f(x_0,y_0)$ exists and define the operator $L\colon X\times Y\longrightarrow Z$ by
\[
L(h,k):=D_1f(x_0,y_0)(h)+\delta_2f(x_0,y_0)(k).
\]
The operator $L$ is linear and continuous. We will show that it is the Fréchet approximation to $f$ at $(x_0,y_0)$.

Let $\varepsilon>0$. By continuity of $\delta_2f$ at $(x_0,y_0)$, there exists $\eta\in(0,r)$ such that
$\|\delta_2f(u,v)-\delta_2f(x_0,y_0)\|_{\mathcal L(Y,Z)}<\displaystyle\frac{\varepsilon}{2}$ whenever $(u,v)\in W$ and $\|(u-x_0,v-y_0)\|_{X\times Y}<\eta$. Since $D_1f(x_0,y_0)$ exists, there exists $\rho>0$ such that
\[
\|f(x_0+h,y_0)-f(x_0,y_0)-D_1f(x_0,y_0)(h)\|_Z
\leq
\frac{\varepsilon}{2}\|h\|_X
\]
whenever $0<\|h\|_X<\rho$; for $h=0$, both sides are zero.

Let $\displaystyle \delta:=\min\{\eta,\rho\}$ and suppose $0<\|h\|_X+\|k\|_Y<\delta$. Then, for every $\tau\in[0,1]$, we have $\|(h,\tau k)\|_{X\times Y}<\eta<r$, so the segment joining $(x_0+h,y_0)$ to $(x_0+h,y_0+k)$ is contained in $W$. Adding and subtracting $f(x_0+h,y_0)$ gives
\[f(x_0+h,y_0+k)-f(x_0,y_0)-L(h,k)\]\[=f(x_0+h,y_0+k)-f(x_0+h,y_0)-\delta_2f(x_0,y_0)(k)\quad+f(x_0+h,y_0)-f(x_0,y_0)-D_1f(x_0,y_0)(h).
\]
Apply the first mean value inequality again, now to the map $v\longmapsto f(x_0+h,v)-\delta_2f(x_0,y_0)(v-y_0)$ between $y_0$ and $y_0+k$. This gives
\[
\begin{aligned}
&\|f(x_0+h,y_0+k)-f(x_0+h,y_0)-\delta_2f(x_0,y_0)(k)\|_Z\\
&\quad\leq\displaystyle\sup_{\tau\in[0,1]}
\|\delta_2f(x_0+h,y_0+\tau k)-\delta_2f(x_0,y_0)\|_{\mathcal L(Y,Z)}\|k\|_Y\\
&\quad\leq\frac{\varepsilon}{2}\|k\|_Y.
\end{aligned}
\]
Consequently,
\[
\|f(x_0+h,y_0+k)-f(x_0,y_0)-L(h,k)\|_Z
\leq
\frac{\varepsilon}{2}\bigl(\|h\|_X+\|k\|_Y\bigr).
\]
It follows that
\[
\lim_{(h,k)\to(0,0)}
\frac{\|f(x_0+h,y_0+k)-f(x_0,y_0)-L(h,k)\|_Z}
{\|h\|_X+\|k\|_Y}
=0.
\]
Thus $f$ is Fréchet differentiable at $(x_0,y_0)$, $Df(x_0,y_0)=L$, and \eqref{eq:descomposicion-derivada-total-frechet} holds.

Finally, Proposition~\ref{prop:derivada-total-implica-parciales} implies that $D_2f(x_0,y_0)$ exists and that, for every $k\in Y$,
\[
D_2f(x_0,y_0)(k)
=
Df(x_0,y_0)(0,k)
=
\delta_2f(x_0,y_0)(k).
\]
Therefore, $D_2f(x_0,y_0)=\delta_2f(x_0,y_0)$.
\end{proof}

\begin{corollary}[Characterization of functions of class $C^1$ by partial derivatives]
\label{cor:C1-derivadas-parciales}
\index{function of class $C^1$!characterization by partial derivatives}
Let $X$, $Y$, and $Z$ be normed spaces, let $U\subseteq X$ and $V\subseteq Y$ be open, and let $f\colon U\times V\longrightarrow Z$. Then $f\in C^1(U\times V,Z)$ if and only if $\delta_1f$ and $\delta_2f$ exist on $U\times V$ and are continuous as maps with values in $\mathcal L(X,Z)$ and $\mathcal L(Y,Z)$, respectively. In this case, $D_1f=\delta_1f$, $D_2f=\delta_2f$, and, for every $(x,y)\in U\times V$ and every $(h,k)\in X\times Y$,
\[
Df(x,y)(h,k)=\delta_1f(x,y)(h)+\delta_2f(x,y)(k).
\]
\end{corollary}

\begin{proof}
First suppose $f\in C^1(U\times V,Z)$. Let $\iota_X\colon X\longrightarrow X\times Y$ and $\iota_Y\colon Y\longrightarrow X\times Y$ be the continuous linear operators defined by $\iota_X(h):=(h,0)$ and $\iota_Y(k):=(0,k)$. By Proposition~\ref{prop:derivada-total-implica-parciales}, for each $(x,y)\in U\times V$ we have $D_1f(x,y)=Df(x,y)\circ\iota_X$ and $D_2f(x,y)=Df(x,y)\circ\iota_Y$. Since $Df\colon U\times V\longrightarrow\mathcal L(X\times Y,Z)$ is continuous, so are $D_1f$ and $D_2f$. Moreover, Fréchet differentiability implies Gâteaux differentiability, and therefore $D_jf=\delta_jf$ for $j\in\{1,2\}$.

Conversely, suppose $\delta_1f$ and $\delta_2f$ exist and are continuous. For fixed $y\in V$, the map $x\longmapsto f(x,y)$ has Gâteaux derivative $\delta_1f(x,y)$, continuous in the operator norm with respect to $x$. Corollary~\ref{cor:caracterizacion-C1-gateaux} shows that this map is of class $C^1$; in particular, $D_1f(x,y)$ exists and agrees with $\delta_1f(x,y)$. All the assumptions of Proposition~\ref{prop:diferenciabilidad-mediante-derivadas-parciales} are now satisfied. Applying it at each point of $U\times V$, we obtain that $f$ is Fréchet differentiable and
\[
Df(x,y)(h,k)=\delta_1f(x,y)(h)+\delta_2f(x,y)(k).
\]
Let $(x_1,y_1),(x_2,y_2)\in U\times V$. For every $(h,k)\in X\times Y$ with $\|(h,k)\|_{X\times Y}\leq1$, we have
\[
\begin{aligned}
\bigl\|\bigl(Df(x_1,y_1)-Df(x_2,y_2)\bigr)(h,k)\bigr\|_Z
&\leq\|\delta_1f(x_1,y_1)-\delta_1f(x_2,y_2)\|_{\mathcal L(X,Z)}\|h\|_X\\
&\quad+\|\delta_2f(x_1,y_1)-\delta_2f(x_2,y_2)\|_{\mathcal L(Y,Z)}\|k\|_Y.
\end{aligned}
\]
Taking the supremum over the unit ball of $X\times Y$ gives
\[
\|Df(x_1,y_1)-Df(x_2,y_2)\|_{\mathcal L(X\times Y,Z)}
\leq
\|\delta_1f(x_1,y_1)-\delta_1f(x_2,y_2)\|_{\mathcal L(X,Z)}
+
\|\delta_2f(x_1,y_1)-\delta_2f(x_2,y_2)\|_{\mathcal L(Y,Z)}.
\]
Continuity of the partial derivatives then implies continuity of $Df$.
\end{proof}

\section{Higher derivatives and Taylor's formula}
\label{sec:derivadas-superiores-banach}

The first derivative takes values in $\mathcal L(X,Y)$. Differentiating it again produces operators that depend linearly on two directions, and further iterations lead to multilinear operators. Before defining higher derivatives, it is therefore useful to specify the normed structure of these spaces and their relation to the spaces of linear operators already introduced.

\begin{definition}[Continuous multilinear operators]
\label{def:operadores-multilineales-continuos}
\index{continuous multilinear operator}
Let $X_1,\ldots,X_m$ and $Y$ be normed spaces, with $m\in\mathbb N$. We denote by $\mathcal L^m(X_1,\ldots,X_m;Y)$ the vector space of continuous $m$-linear maps $A\colon X_1\times\cdots\times X_m\longrightarrow Y$. For each $A\in\mathcal L^m(X_1,\ldots,X_m;Y)$, define
\[
\|A\|_{\mathcal L^m}
:=
\sup\bigl\{\|A(h_1,\ldots,h_m)\|_Y\mid h_j\in X_j,\ \|h_j\|_{X_j}\leq1\text{ for }1\leq j\leq m\bigr\}.
\]
If $X_1=\cdots=X_m=X$, write $\mathcal L^m(X;Y)$. In particular, $\mathcal L^1(X;Y)=\mathcal L(X,Y)$. Denote by $\mathcal L_s^m(X;Y)$ the subspace of symmetric operators, that is, those $A\in\mathcal L^m(X;Y)$ such that $A(h_{\sigma(1)},\ldots,h_{\sigma(m)})=A(h_1,\ldots,h_m)$ for every permutation $\sigma$ of $\{1,\ldots,m\}$. By convention, $\mathcal L^0(X;Y):=Y$.
\end{definition}

\begin{proposition}[Characterization of multilinear continuity]
\label{prop:caracterizacion-multilineales-continuos}
A $m$-linear map $A\colon X_1\times\cdots\times X_m\longrightarrow Y$ is continuous if and only if there exists a constant $C\geq0$ such that
\begin{equation}
\label{eq:cota-operador-multilineal}
\|A(h_1,\ldots,h_m)\|_Y
\leq
C\prod_{j=1}^{m}\|h_j\|_{X_j}
\end{equation}
for all $h_j\in X_j$. In this case, the smallest constant for which \eqref{eq:cota-operador-multilineal} holds is $\|A\|_{\mathcal L^m}$.
\end{proposition}

\begin{proof}
First suppose \eqref{eq:cota-operador-multilineal} holds. If $h_j\longrightarrow0$ in $X_j$ for each $j$, the right-hand side of \eqref{eq:cota-operador-multilineal} tends to zero, so $A$ is continuous at the origin. Now let $x_j,h_j\in X_j$. By multilinearity,
\[
A(x_1+h_1,\ldots,x_m+h_m)-A(x_1,\ldots,x_m)
=
\displaystyle\sum_{\varnothing\neq I\subseteq\{1,\ldots,m\}}A(z_1^I,\ldots,z_m^I),
\]
where $z_j^I=h_j$ if $j\in I$ and $z_j^I=x_j$ if $j\notin I$. Inequality \eqref{eq:cota-operador-multilineal}, applied to each summand, shows that the difference tends to zero as $(h_1,\ldots,h_m)\longrightarrow(0,\ldots,0)$. Thus $A$ is continuous at every point.

Conversely, suppose $A$ is continuous at the origin. There exists $\rho>0$ such that $\|A(u_1,\ldots,u_m)\|_Y<1$ whenever $\|u_j\|_{X_j}<\rho$ for every $j$. If none of the vectors $h_j$ is zero, define $\displaystyle u_j:=\displaystyle\frac{\rho h_j}{2\|h_j\|_{X_j}}$. Then $\|u_j\|_{X_j}=\displaystyle\frac{\rho}{2}$ and, by multilinearity,
\[
\|A(h_1,\ldots,h_m)\|_Y
<
\left(\frac{2}{\rho}\right)^m\prod_{j=1}^{m}\|h_j\|_{X_j}.
\]
If one of the vectors is zero, both sides of \eqref{eq:cota-operador-multilineal} vanish. This proves existence of a constant $C$.

Finally, if \eqref{eq:cota-operador-multilineal} holds with constant $C$, then $\|A\|_{\mathcal L^m}\leq C$. On the other hand, for nonzero vectors, applying the definition of the norm to $\displaystyle\frac{h_j}{\|h_j\|_{X_j}}$ gives
\[
\|A(h_1,\ldots,h_m)\|_Y
\leq
\|A\|_{\mathcal L^m}\prod_{j=1}^{m}\|h_j\|_{X_j}.
\]
Thus $\|A\|_{\mathcal L^m}$ is the smallest possible constant.
\end{proof}

\begin{proposition}[Completeness and association of multilinear operators]
\label{prop:completitud-multilineales}
If $Y$ is Banach, then $\mathcal L^m(X_1,\ldots,X_m;Y)$ is Banach. Moreover, the map
\[
\mathfrak C:\mathcal L^{m+1}(X_1,\ldots,X_m,X_{m+1};Y)
\longrightarrow
\mathcal L\bigl(X_{m+1},\mathcal L^m(X_1,\ldots,X_m;Y)\bigr),
\]
defined by
\[
\mathfrak C(A)(h_{m+1})(h_1,\ldots,h_m)
:=
A(h_1,\ldots,h_m,h_{m+1}),
\]
is an isometric linear isomorphism. In particular, if $Y$ is Banach, then $\mathcal L_s^m(X;Y)$ is a closed subspace and hence a Banach space.
\end{proposition}

\begin{proof}
Let $(A_j)_{j\in\mathbb N}$ be a Cauchy sequence in $\mathcal L^m(X_1,\ldots,X_m;Y)$. For each $(h_1,\ldots,h_m)\in X_1\times\cdots\times X_m$, we have
\[
\|(A_j-A_k)(h_1,\ldots,h_m)\|_Y
\leq
\|A_j-A_k\|_{\mathcal L^m}\prod_{i=1}^{m}\|h_i\|_{X_i}.
\]
Thus $(A_j(h_1,\ldots,h_m))_{j\in\mathbb N}$ is a Cauchy sequence in $Y$. Define
$\displaystyle A(h_1,\ldots,h_m):=\lim_{j\to\infty}A_j(h_1,\ldots,h_m)$. Passage to the limit in each variable shows that $A$ is $m$-linear.

Since $(A_j)$ is Cauchy, there exists $j_0\in\mathbb N$ such that $\|A_j-A_k\|_{\mathcal L^m}<1$ for all $j,k\geq j_0$. Letting $k\to\infty$ in the preceding estimate gives
\[
\|(A_j-A)(h_1,\ldots,h_m)\|_Y
\leq
\prod_{i=1}^{m}\|h_i\|_{X_i}
\]
for every $j\geq j_0$. Proposition~\ref{prop:caracterizacion-multilineales-continuos} implies that $A_j-A$ is continuous. Moreover, given $\varepsilon>0$, there exists $j_1\in\mathbb N$ such that $\|A_j-A_k\|_{\mathcal L^m}<\varepsilon$ for all $j,k\geq j_1$. Letting $k\to\infty$ and taking the supremum over the unit balls gives $\|A_j-A\|_{\mathcal L^m}\leq\varepsilon$ for every $j\geq j_1$. Thus $A_j\longrightarrow A$, and the space is complete.

It is immediate that $\mathfrak C$ is linear. For each $A$, we have
\[
\|\mathfrak C(A)\|
=
\sup_{\substack{h_{m+1}\in X_{m+1}\\\|h_{m+1}\|_{X_{m+1}}\leq1}}\sup_{\substack{h_i\in X_i,\ 1\leq i\leq m\\\|h_i\|_{X_i}\leq1}}
\|A(h_1,\ldots,h_m,h_{m+1})\|_Y
=
\|A\|_{\mathcal L^{m+1}},
\]
so it is an isometry. If $T\in\mathcal L(X_{m+1},\mathcal L^m(X_1,\ldots,X_m;Y))$, the map
$A_T(h_1,\ldots,h_{m+1}):=T(h_{m+1})(h_1,\ldots,h_m)$ is continuous and $(m+1)$-linear and satisfies $\mathfrak C(A_T)=T$. Thus $\mathfrak C$ is surjective.

Finally, if a sequence of symmetric operators converges in norm to $A$, passage to the limit in the symmetry identity shows that $A$ is also symmetric. Thus $\mathcal L_s^m(X;Y)$ is closed.
\end{proof}

For $A\in\mathcal L^m(X;Y)$, write $A[h_1,\ldots,h_m]:=A(h_1,\ldots,h_m)$ and $A[h]^m:=A(h,\ldots,h)$.

\begin{proposition}[Differentiability of multilinear operators]
\label{prop:diferenciabilidad-operadores-multilineales}
Let $X_1,\ldots,X_m$ and $Y$ be normed spaces, and let $A\in\mathcal L^m(X_1,\ldots,X_m;Y)$. Regarded as a map on the product space, $A$ is smooth and
\[
DA(x_1,\ldots,x_m)(h_1,\ldots,h_m)
=
\displaystyle\sum_{j=1}^{m}A(x_1,\ldots,x_{j-1},h_j,x_{j+1},\ldots,x_m).
\]
\end{proposition}

\begin{proof}
By multilinearity,
\[
A(x_1+h_1,\ldots,x_m+h_m)-A(x_1,\ldots,x_m)
\]
is the sum of the terms in which at least one of the variables $x_j$ has been replaced by $h_j$. The terms containing exactly one increment form the linear map in the statement. The remaining terms contain at least two increments. If $\displaystyle \max_{j\in\{1,\ldots,m\}}\|h_j\|_{X_j}<1$, Proposition~\ref{prop:caracterizacion-multilineales-continuos} bounds the norm of their sum by
\[
\|A\|_{\mathcal L^m}
\displaystyle\sum_{\substack{I\subseteq\{1,\ldots,m\}\\\#I\geq2}}
\prod_{j\in I}\|h_j\|_{X_j}
\prod_{j\in\{1,\ldots,m\}\setminus I}\max\{1,\|x_j\|_{X_j}\}.
\]
After division by $\displaystyle\sum_{j=1}^{m}\|h_j\|_{X_j}$, each summand tends to zero as $(h_1,\ldots,h_m)\longrightarrow0$. This proves the formula for the first derivative. To make the successive derivatives explicit, write $E:=X_1\times\cdots\times X_m$, with the sum norm, and fix $k\in\{1,\ldots,m\}$. Let $\mathcal I_k$ be the set of injective maps $\sigma\colon\{1,\ldots,k\}\longrightarrow\{1,\ldots,m\}$. For $x=(x_1,\ldots,x_m)\in E$ and $h^{(a)}=(h_1^{(a)},\ldots,h_m^{(a)})\in E$, with $a\in\{1,\ldots,k\}$, set
\[
z_i^\sigma:=
\begin{cases}
h_i^{(a)},&i=\sigma(a)\text{ for some }a\in\{1,\ldots,k\},\\
x_i,&i\in\{1,\ldots,m\}\setminus\operatorname{im}(\sigma).
\end{cases}
\]
Injectivity of $\sigma$ makes this definition unambiguous. Define the continuous $k$--linear form
\[
B_k(x)[h^{(1)},\ldots,h^{(k)}]
:=\sum_{\sigma\in\mathcal I_k}A(z_1^\sigma,\ldots,z_m^\sigma).
\]
The formula already proved says that $DA=B_1$. Now suppose $1\leq k<m$. To compute the increment of $B_k$ when $x$ is replaced by $x+u$, expand each summand by multilinearity. The terms linear in $u$ are obtained by choosing an index $i\in\{1,\ldots,m\}\setminus\operatorname{im}(\sigma)$ and replacing $x_i$ by $u_i$. To each pair $(\sigma,i)$ there corresponds exactly one map in $\mathcal I_{k+1}$, agreeing with $\sigma$ on $\{1,\ldots,k\}$ and sending $k+1$ to $i$. The sum of these terms is therefore $B_{k+1}(x)[h^{(1)},\ldots,h^{(k)},u]$.

The remaining terms contain at least two components of $u$. If $\|u\|_E\leq1$, the multilinear bound estimates them together by
\[
C_{A,x,k}\|u\|_E^2\prod_{a=1}^{k}\|h^{(a)}\|_E,
\]
where $C_{A,x,k}$ is finite: only finitely many summands, the norm of $A$, and the factors $\displaystyle \max\{1,\|x_i\|_{X_i}\}$, with $i\in\{1,\ldots,m\}$, occur. Taking the supremum over $\|h^{(a)}\|_E\leq1$ for $a\in\{1,\ldots,k\}$ and dividing by $\|u\|_E$, we obtain that $B_k$ is Fréchet differentiable as a map with values in $\mathcal L^k(E;Y)$. Its derivative corresponds to $B_{k+1}$ through Proposition~\ref{prop:completitud-multilineales}. The same expansion, retaining all terms with at least one increment, proves norm continuity of each $B_k$.

This completes the step from order $k$ to order $k+1$. For $k=m$, all indices have been occupied by the increments $h^{(a)}$, so $B_m$ is independent of $x$ and its derivative is zero. Subsequent derivatives are also zero. Consequently, $A$ is smooth.
\end{proof}

\begin{definition}[Higher Gâteaux derivatives]
\label{def:derivadas-superiores-gateaux}
\index{Gateaux derivative@Gâteaux derivative!higher order}
Let $U\subseteq X$ be open and let $f\colon U\longrightarrow Y$. Define $\delta^0f:=f$. Suppose $\delta^{m-1}f$ is defined on a neighborhood of $x\in U$ as a map with values in $\mathcal L^{m-1}(X;Y)$ and is Gâteaux differentiable at $x$. Through the isomorphism in Proposition~\ref{prop:completitud-multilineales}, identify
\[
\delta(\delta^{m-1}f)(x)\in\mathcal L\bigl(X,\mathcal L^{m-1}(X;Y)\bigr)
\]
with an element $\delta^mf(x)\in\mathcal L^m(X;Y)$, called the \textbf{Gâteaux derivative of order $m$} of $f$ at $x$. With the order of association fixed by $\mathfrak C$, we have
\[
\delta^mf(x)[h_1,\ldots,h_m]
=
D_{h_m}\bigl(D_{h_{m-1}}\cdots D_{h_1}f\bigr)(x)
\]
whenever the iterated derivatives on the right-hand side are defined on a neighborhood of $x$.
\end{definition}

\begin{definition}[Higher Fréchet derivatives and functions of class $C^k$]
\label{def:funcion-Ck-banach}
\index{Frechet derivative@Fréchet derivative!higher order}
\index{function of class $C^k$!between Banach spaces}
Let $U\subseteq X$ be open and let $f\colon U\longrightarrow Y$. Define $D^0f:=f$. Suppose $D^{m-1}f$ is defined on a neighborhood of $x\in U$ as a map with values in $\mathcal L^{m-1}(X;Y)$ and is Fréchet differentiable at $x$. Identify $D(D^{m-1}f)(x)$ with an element $D^mf(x)\in\mathcal L^m(X;Y)$ through Proposition~\ref{prop:completitud-multilineales}; it is called the \textbf{Fréchet derivative of order $m$} of $f$ at $x$.

For $k\in\mathbb N$, we say that $f$ is of \textbf{class $C^k$} on $U$, and write $f\in C^k(U,Y)$, if $D^jf$ exists at every point of $U$ and is continuous as a map $U\longrightarrow\mathcal L^j(X;Y)$ for every $j\in\{1,\ldots,k\}$. We say that $f$ is \textbf{smooth}, and write $f\in C^\infty(U,Y)$, if $f\in C^k(U,Y)$ for every $k\in\mathbb N$.
\end{definition}

\begin{remark}
In the definition, continuity of all derivatives through order $k$ can be replaced by existence of $D^kf$ and continuity of the latter. For each $j<k$, the map $D^jf$ is Fréchet differentiable and hence continuous by Proposition~\ref{prop:frechet-implica-continuidad}. Nevertheless, the preceding formulation makes all the available regularity explicit and is more convenient when applying higher-order chain rules.
\end{remark}

Regularity of a composition follows from the first-order chain rule and continuity of operator composition. It is useful to formulate this fact with values in normed spaces, since the derivatives themselves take values in operator spaces.

\begin{corollary}[Composition of maps of class $C^k$]
\label{cor:composicion-Ck-banach}
Let $X$, $Y$, and $Z$ be normed spaces, let $U\subseteq X$ and $V\subseteq Y$ be open, and let $k\in\mathbb N\cup\{\infty\}$. If $f\in C^k(U,Y)$, $f(U)\subseteq V$, and $g\in C^k(V,Z)$, then $g\circ f\in C^k(U,Z)$.
\end{corollary}

\begin{proof}
The first-order chain rule gives
\[
D(g\circ f)(x)=Dg(f(x))\circ Df(x),
\qquad x\in U.
\]
Operator composition
\[
\mathcal C\colon\mathcal L(Y,Z)\times\mathcal L(X,Y)\longrightarrow\mathcal L(X,Z),
\qquad \mathcal C(B,A):=B\circ A,
\]
is continuous and bilinear, since $\|B\circ A\|\leq\|B\|\|A\|$. If $f$ and $g$ are of class $C^1$, the maps $Dg\circ f$ and $Df$ are continuous. The preceding identity then shows that $D(g\circ f)$ is continuous, proving the case $k=1$.

Suppose the assertion has been proved for an integer $r\geq1$ in arbitrary normed spaces, and now let $f$ and $g$ be of class $C^{r+1}$. The map $Dg$ is of class $C^r$; by the induction hypothesis, so is $Dg\circ f$. Moreover, $Df$ is of class $C^r$. The map $x\longmapsto(Dg(f(x)),Df(x))$ is therefore of class $C^r$: in a product of two normed spaces, differentiation and continuity are checked componentwise using the sum norm. By Proposition~\ref{prop:diferenciabilidad-operadores-multilineales}, $\mathcal C$ is smooth. Applying the induction hypothesis again to composition with $\mathcal C$, we conclude that $D(g\circ f)$ is of class $C^r$. Thus $g\circ f$ is of class $C^{r+1}$. Induction proves all finite orders; the case $k=\infty$ follows by applying the conclusion to each of them.
\end{proof}

\begin{theorem}[Characterization by Gâteaux derivatives]
\label{teo:caracterizacion-Ck-gateaux}
\index{function of class $C^k$!characterization by Gâteaux derivatives}
Let $X$ and $Y$ be normed spaces, let $U\subseteq X$ be open, and let $k\in\mathbb N$. Then $f\in C^k(U,Y)$ if and only if, for each $j\in\{1,\ldots,k\}$, the Gâteaux derivative $\delta^jf(x)\in\mathcal L^j(X;Y)$ exists for every $x\in U$ and the map $\delta^jf\colon U\longrightarrow\mathcal L^j(X;Y)$ is continuous. In this case, $D^jf=\delta^jf$ for every $j\in\{1,\ldots,k\}$.
\end{theorem}

\begin{proof}
We proceed by induction on $k$. For $k=1$, the assertion is Corollary~\ref{cor:caracterizacion-C1-gateaux}. Suppose the result holds for $k-1$.

If $f\in C^k(U,Y)$, then $D^{k-1}f\colon U\longrightarrow\mathcal L^{k-1}(X;Y)$ is of class $C^1$. Corollary~\ref{cor:caracterizacion-C1-gateaux}, applied to this map, shows that its Gâteaux derivative exists, is continuous, and agrees with its Fréchet derivative. Through the identification in Proposition~\ref{prop:completitud-multilineales}, this means that $\delta^kf=D^kf$. For lower orders, the conclusion follows from the induction hypothesis.

Conversely, suppose the Gâteaux derivatives in the statement exist and are continuous. By the induction hypothesis, $f\in C^{k-1}(U,Y)$ and $D^{k-1}f=\delta^{k-1}f$. The map $D^{k-1}f$ is Gâteaux differentiable, its derivative is $\delta^kf$, and this derivative is continuous. Corollary~\ref{cor:caracterizacion-C1-gateaux}, applied with target space $\mathcal L^{k-1}(X;Y)$, implies that $D^{k-1}f$ is of class $C^1$ and its Fréchet derivative agrees with $\delta^kf$. Thus $f\in C^k(U,Y)$ and $D^kf=\delta^kf$.
\end{proof}

\begin{theorem}[Symmetry of higher derivatives]
\label{teo:simetria-derivadas-superiores}
\index{Frechet derivative@Fréchet derivative!symmetry}
If $f\in C^m(U,Y)$, with $m\geq2$, then $D^mf(x)\in\mathcal L_s^m(X;Y)$ for every $x\in U$.
\end{theorem}

\begin{proof}
We proceed by induction on $m$. For $m=2$, fix $x\in U$, $h,k\in X$, and $\varphi\in Y'$. Define $g(s,t):=\varphi(f(x+sh+tk))$ for $(s,t)$ in a neighborhood of $(0,0)$. The chain rule implies that $g$ is of class $C^2$ and that
\[
\frac{\partial^2g}{\partial s\,\partial t}(0,0)
=
\varphi(D^2f(x)[k,h]),
\qquad
\frac{\partial^2g}{\partial t\,\partial s}(0,0)
=
\varphi(D^2f(x)[h,k]).
\]
The classical equality of mixed partial derivatives and the fact that $Y'$ separates points imply $D^2f(x)[h,k]=D^2f(x)[k,h]$.

Now suppose the assertion holds at order $m-1$, with $m\geq3$. Since $D^{m-1}f(x)$ is symmetric, the operator $D^mf(x)$ is symmetric in its first $m-1$ arguments. For each $v\in X$, the function
\[
x\longmapsto D^{m-1}f(x)[h_{\sigma(1)},\ldots,h_{\sigma(m-1)}]
-
D^{m-1}f(x)[h_1,\ldots,h_{m-1}]
\]
is identically zero for every permutation $\sigma$ of $\{1,\ldots,m-1\}$; differentiating it in the direction $v$ gives the stated symmetry.

It remains to show that the last two arguments can be interchanged. Fix $h_1,\ldots,h_{m-2}\in X$ and consider
\[
G(z):=D^{m-2}f(z)[h_1,\ldots,h_{m-2}].
\]
Evaluation of a multilinear operator on fixed vectors is a continuous linear map; hence $G$ is of class $C^2$. The case $m=2$, applied to $G$, gives
\[
D^m f(x)[h_1,\ldots,h_{m-2},h_{m-1},h_m]
=
D^m f(x)[h_1,\ldots,h_{m-2},h_m,h_{m-1}].
\]
Permutations of the first $m-1$ arguments, together with transposition of the last two, generate all permutations of $m$ elements. Thus $D^mf(x)$ is symmetric.
\end{proof}

\begin{corollary}[Faà di Bruno for the Fréchet derivative]
\label{cor:faa-di-bruno-espacios-normados}
Let $X,Y,Z$ be normed spaces, let $U\subseteq X$ and
$V\subseteq Y$ be open, and let $k\geq1$ be an integer.
Let $f\in C^k(U,Y)$ and $g\in C^k(V,Z)$ satisfy $f(U)\subseteq V$.
In particular, $X,Y,Z$ may be Banach spaces. All derivatives
below are Fréchet derivatives.

For $1\leq q\leq k$, let $\Pi_q$ be the set of partitions of
$\{1,\ldots,q\}$ into nonempty blocks. Enumerate the blocks of
$\pi=\{B_1,\ldots,B_a\}$ so that $\min B_1<\cdots<\min B_a$.
If $v=(v_1,\ldots,v_q)\in X^q$ and $B=\{i_1<\cdots<i_b\}$, write
\[
 D_Bf(x)[v]:=D^bf(x)[v_{i_1},\ldots,v_{i_b}].
\]
Then, for every $x\in U$,
\[
 D^q(g\circ f)(x)[v_1,\ldots,v_q]
 =\sum_{\substack{\pi\in\Pi_q\\\pi=\{B_1,\ldots,B_a\}}}
 D^ag(f(x))\bigl[D_{B_1}f(x)[v],\ldots,D_{B_a}f(x)[v]\bigr].
\]
This identity is an equality of operators in $\mathcal L^q(X;Z)$
evaluated at $v_1,\ldots,v_q$.
If $I\subseteq\mathbb R$ is an open interval and
$f\in C^k(I,Y)$ satisfies $f(I)\subseteq V$, then, for $t\in I$,
\[
 (g\circ f)^{(q)}(t)
 =\sum_{j=1}^{q}
 \sum_{\substack{k_1,\ldots,k_j\geq1\\k_1+\cdots+k_j=q}}
 \frac{q!}{j!\,k_1!\cdots k_j!}
 D^jg(f(t))\bigl[f^{(k_1)}(t),\ldots,f^{(k_j)}(t)\bigr],
\]
where the indices $k_1,\ldots,k_j$ are integers.
For smooth maps the identities hold for every integer $q\geq1$.
\end{corollary}

\begin{proof}
By Corollary~\ref{cor:composicion-Ck-banach}, the composition
is of class $C^k$. Its derivatives and those of $f$ and $g$
are symmetric by Theorem~\ref{teo:simetria-derivadas-superiores}.
We first justify the required product rule as an identity
in operator norm.

Fix $q\geq1$ and a partition
$\pi=\{B_1,\ldots,B_a\}\in\Pi_q$, and put $b_\ell=|B_\ell|$.
For $T\in\mathcal L^{b_\ell}(X;Y)$ and
$B_\ell=\{i_1<\cdots<i_{b_\ell}\}$, use the abbreviation
$T[v_{B_\ell}]:=T[v_{i_1},\ldots,v_{i_{b_\ell}}]$.
Define the map
\[
\begin{split}
 \mathcal C_\pi\colon
 \mathcal L^a(Y;Z)\times
 \prod_{\ell=1}^{a}\mathcal L^{b_\ell}(X;Y)
 &\longrightarrow\mathcal L^q(X;Z),\\
 \mathcal C_\pi(A,T_1,\ldots,T_a)[v_1,\ldots,v_q]
 &:=A[T_1[v_{B_1}],\ldots,T_a[v_{B_a}]].
\end{split}
\]
Each direction $v_i$ occurs in exactly one block. The expression
is therefore $q$-linear in the directions. Moreover,
\[
 \|\mathcal C_\pi(A,T_1,\ldots,T_a)[v_1,\ldots,v_q]\|_Z
 \leq\|A\|_{\mathcal L^a(Y;Z)}
       \prod_{\ell=1}^{a}\|T_\ell\|_{\mathcal L^{b_\ell}(X;Y)}
       \prod_{i=1}^{q}\|v_i\|_X.
\]
Taking the supremum over $\|v_i\|_X\leq1$ gives
\[
 \|\mathcal C_\pi(A,T_1,\ldots,T_a)\|_{\mathcal L^q(X;Z)}
 \leq\|A\|_{\mathcal L^a(Y;Z)}
       \prod_{\ell=1}^{a}\|T_\ell\|_{\mathcal L^{b_\ell}(X;Y)}.
\]
Thus $\mathcal C_\pi$ is well defined and is a continuous
$(a+1)$-linear map. Proposition~\ref{prop:diferenciabilidad-operadores-multilineales}
applies to this map: its derivative is the sum of the terms
obtained by incrementing one argument at a time. In particular,
this rule holds in the norm of $\mathcal L^q(X;Z)$.

We now prove the formula by induction. For $q=1$, it is the
chain rule. Suppose it holds for $q<k$. For each partition
$\pi$ as above, define
\[
 A(x):=D^ag(f(x)),\qquad T_\ell(x):=D^{b_\ell}f(x),\qquad
 S_\pi(x):=\mathcal C_\pi(A(x),T_1(x),\ldots,T_a(x)).
\]
Since $a,b_1,\ldots,b_a\leq q<k$, these factors are of class
$C^1$ in their respective operator spaces. For $h\in X$ and
$y_1,\ldots,y_a\in Y$, the chain rule and the identification
of iterated derivatives with multilinear operators give
\[
\begin{split}
 (DA(x)[h])[y_1,\ldots,y_a]
   &=D^{a+1}g(f(x))[y_1,\ldots,y_a,Df(x)[h]],\\
 (DT_\ell(x)[h])[v_{B_\ell}]
   &=D^{b_\ell+1}f(x)[v_{i_1},\ldots,v_{i_{b_\ell}},h].
\end{split}
\]
These are Fréchet derivatives in operator norm, by the definition
of $C^k$ and the isometric identification used in
Proposition~\ref{prop:completitud-multilineales}.
The chain rule applied to $\mathcal C_\pi$ shows that
$S_\pi\colon U\to\mathcal L^q(X;Z)$ is of class $C^1$, with
\[
\begin{split}
 DS_\pi(x)[h]
 ={}&\mathcal C_\pi(DA(x)[h],T_1(x),\ldots,T_a(x))\\
 &+\sum_{\ell=1}^{a}\mathcal C_\pi
 (A(x),T_1(x),\ldots,T_{\ell-1}(x),DT_\ell(x)[h],
                  T_{\ell+1}(x),\ldots,T_a(x)).
\end{split}
\]
To describe the terms obtained by evaluating this identity, write
\[
 W_\ell:=D_{B_\ell}f(x)[v],\qquad
 W_\ell^h:=D^{b_\ell+1}f(x)
             [v_{i_1},\ldots,v_{i_{b_\ell}},h].
\]
We obtain
\[
\begin{split}
 (DS_\pi(x)[h])[v_1,\ldots,v_q]
 ={}&D^{a+1}g(f(x))[W_1,\ldots,W_a,Df(x)[h]]\\
 &+\sum_{\ell=1}^{a}D^ag(f(x))
       [W_1,\ldots,W_{\ell-1},W_\ell^h,W_{\ell+1},\ldots,W_a].
\end{split}
\]
The induction hypothesis is the equality
$D^q(g\circ f)(x)=\displaystyle\sum_{\pi\in\Pi_q}S_\pi(x)$
in $\mathcal L^q(X;Z)$ for every $x\in U$. Since the sum is
finite, we may differentiate it in that space. Evaluating the
derivative at $h=v_{q+1}$ and then at $v_1,\ldots,v_q$ gives
$D^{q+1}(g\circ f)(x)[v_1,\ldots,v_q,v_{q+1}]$ on the left.

On the right, the first term of each $DS_\pi$ corresponds to
adjoining the block $\{q+1\}$. The term with index $\ell$
corresponds to adjoining $q+1$ to the block $B_\ell$.
To check that every partition of $\{1,\ldots,q+1\}$ occurs
exactly once, take such a partition and consider its block $C$
containing $q+1$. If $C$ is a singleton, deleting it recovers
the unique original partition. If $|C|\geq2$, replacing it by
$C\setminus\{q+1\}$ recovers both the original partition and
its distinguished block. The alternatives are disjoint and
exhaust all cases. Symmetry of the derivatives allows the blocks
to be ordered as in the statement. This proves the formula for
$q+1$ and completes the induction.

Finally, for a curve evaluate every direction at $1\in\mathbb R$.
Each block of size $k_\ell$ contributes $f^{(k_\ell)}(t)$.
For $j$ blocks with ordered sizes $k_1,\ldots,k_j\geq1$ whose
sum is $q$, the number of successive choices of their positions is
\[
 \binom q{k_1}\binom{q-k_1}{k_2}\cdots
 \binom{q-k_1-\cdots-k_{j-1}}{k_j}
 =\frac{q!}{k_1!\cdots k_j!}.
\]
Each unordered partition occurs $j!$ times among all ordered
lists of blocks. This remains true when some sizes coincide,
because the blocks are still distinct subsets. Symmetry of
$D^jg(f(t))$ permits division by $j!$. Summing over the sizes
and over $j\in\{1,\ldots,q\}$ gives the curve formula.
The proof does not use completeness, so it establishes the
result for normed spaces as well as for Banach spaces.
\end{proof}

The Leibniz rule has a particularly useful formulation for the second derivative. We will use this identity in the next chapter when comparing the kinematic and operational notions of tangent space.

\begin{proposition}[Leibniz rule for the second derivative]
\label{prop:leibniz-segunda-derivada}
\index{Leibniz rule!second derivative}
Let $X$ be a normed space, let $U\subseteq X$ be open, and let $f,g\in C^2(U,\mathbb R)$. Then, for every $x\in U$ and $h,k\in X$,
\[
D^2(fg)(x)[h,k]
=
f(x)D^2g(x)[h,k]+g(x)D^2f(x)[h,k]
+Df(x)(h)Dg(x)(k)+Df(x)(k)Dg(x)(h).
\]
Equivalently,
\[
D^2(fg)(x)
=
f(x)D^2g(x)+g(x)D^2f(x)+Df(x)\otimes Dg(x)+Dg(x)\otimes Df(x).
\]
\end{proposition}

\begin{proof}
The product rule for the first derivative gives $D(fg)(x)=f(x)Dg(x)+g(x)Df(x)$. Differentiate this map with values in $X'$ in the direction $k$ and evaluate the resulting operator on $h$. The product rule and linearity of evaluation give
\[D^2(fg)(x)[h,k]
=Df(x)(k)Dg(x)(h)+f(x)D^2g(x)[h,k]+Dg(x)(k)Df(x)(h)+g(x)D^2f(x)[h,k].
\]
Symmetry of second derivatives, given by Theorem~\ref{teo:simetria-derivadas-superiores}, allows the arguments to be ordered as in the statement.
\end{proof}

\begin{theorem}[Taylor's formula with integral remainder]
\label{teo:taylor-banach}
\index{Taylor's formula!in Banach spaces}
Let $X$ be a normed space, let $Y$ be a Banach space, let $U\subseteq X$ be open, and let $m\in\mathbb N$ and $f\in C^m(U,Y)$. Let $x\in U$ and $h\in X$ be such that $[x,x+h]\subseteq U$. Then
\begin{equation}
\label{eq:taylor-resto-integral-banach}
f(x+h)
=
\displaystyle\sum_{j=0}^{m-1}\frac{1}{j!}D^jf(x)[h]^j
+
\frac{1}{(m-1)!}\int_0^1(1-t)^{m-1}D^mf(x+th)[h]^m\,dt,
\end{equation},
where $D^0f(x)[h]^0:=f(x)$.
\end{theorem}

\begin{proof}
Define $\varphi\colon [0,1]\longrightarrow Y$ by $\varphi(t):=f(x+th)$. Since the segment $[x,x+h]$ is contained in the open set $U$, this formula defines the curve on an open interval containing $[0,1]$. The chain rule gives $\varphi'(t)=Df(x+th)(h)$. More generally, we will prove that $\varphi^{(j)}(t)=D^jf(x+th)[h]^j$ for every $j\in\{0,\ldots,m\}$.

The case $j=0$ is the definition of $\varphi$, and the case $j=1$ has just been checked. Suppose the formula holds for an integer $j$ with $1\leq j<m$. The evaluation $A\longmapsto A[h]^j$ is a continuous linear operator from $\mathcal L^j(X;Y)$ to $Y$, with norm at most $\|h\|_X^j$. We may therefore differentiate the inductive identity by the chain rule. The derivative of $D^jf(x+th)$ with respect to $t$ is $D(D^jf)(x+th)(h)$. Evaluating on the $j$ vectors $h$ and using the association in Proposition~\ref{prop:completitud-multilineales} gives $D^{j+1}f(x+th)[h]^{j+1}$. This proves the induction step. Continuity of these expressions follows from continuity of each $D^jf$.

We prove by induction that every function $\varphi\in C^m([0,1],Y)$ satisfies
\begin{equation}
\label{eq:taylor-curva-banach}
\varphi(1)
=
\displaystyle\sum_{j=0}^{m-1}\frac{\varphi^{(j)}(0)}{j!}
+
\frac{1}{(m-1)!}\int_0^1(1-t)^{m-1}\varphi^{(m)}(t)\,dt.
\end{equation}
For $m=1$, this identity is Theorem~\ref{teo:b5-fundamental-calculo-riemann-banach}. Suppose it holds for $m-1$. The induction hypothesis gives
\[
\varphi(1)
=
\displaystyle\sum_{j=0}^{m-2}\frac{\varphi^{(j)}(0)}{j!}
+
\frac{1}{(m-2)!}\int_0^1(1-t)^{m-2}\varphi^{(m-1)}(t)\,dt.
\]
Apply the integration by parts formula in Proposition~\ref{prop:b5-integracion-por-partes-riemann-banach} to the scalar function $\displaystyle a(t):=-\displaystyle\frac{(1-t)^{m-1}}{(m-1)!}$ and the vector-valued function $\varphi^{(m-1)}$. Since $\displaystyle a'(t)=\displaystyle\frac{(1-t)^{m-2}}{(m-2)!}$, we obtain
\[
\frac{1}{(m-2)!}\int_0^1(1-t)^{m-2}\varphi^{(m-1)}(t)\,dt
=
\frac{\varphi^{(m-1)}(0)}{(m-1)!}
+
\frac{1}{(m-1)!}\int_0^1(1-t)^{m-1}\varphi^{(m)}(t)\,dt.
\]
Substituting this equality into the induction hypothesis gives \eqref{eq:taylor-curva-banach}. Finally, replacing each $\varphi^{(j)}(t)$ by $D^jf(x+th)[h]^j$ yields \eqref{eq:taylor-resto-integral-banach}.
\end{proof}

\begin{corollary}[Estimate for the Taylor remainder]
\label{cor:estimacion-resto-taylor-banach}
Under the assumptions of Theorem~\ref{teo:taylor-banach}, suppose $r>0$ satisfies $B_{d_X}(x,r)\subseteq U$. Then, for each $h\in X$ with $\|h\|_X<r$,
\[
\left\|f(x+h)-\displaystyle\sum_{j=0}^{m}\frac{1}{j!}D^jf(x)[h]^j\right\|_Y
\leq
\frac{\|h\|_X^m}{m!}
\sup_{t\in[0,1]}\|D^mf(x+th)-D^mf(x)\|_{\mathcal L^m}.
\]
In particular, for each $\varepsilon>0$ there exists $\delta\in(0,r)$ such that
\[
\left\|f(x+h)-\displaystyle\sum_{j=0}^{m}\frac{1}{j!}D^jf(x)[h]^j\right\|_Y
\leq
\varepsilon\|h\|_X^m
\]
whenever $\|h\|_X<\delta$.
\end{corollary}

\begin{proof}
In \eqref{eq:taylor-resto-integral-banach}, add and subtract $D^mf(x)[h]^m$ inside the integral. Since $\displaystyle\int_0^1(1-t)^{m-1}\,dt=\displaystyle\frac{1}{m}$, we obtain
\[
f(x+h)-\displaystyle\sum_{j=0}^{m}\frac{1}{j!}D^jf(x)[h]^j=
\frac{1}{(m-1)!}\int_0^1(1-t)^{m-1}
\bigl(D^mf(x+th)-D^mf(x)\bigr)[h]^m\,dt.
\]
The inequality for the Riemann integral with values in $Y$ and the characterization of the multilinear norm imply the first estimate. Now let $\varepsilon>0$. By continuity of $D^mf$ at $x$, there exists $\delta\in(0,r)$ such that
\[
\|D^mf(x+u)-D^mf(x)\|_{\mathcal L^m}<m!\varepsilon
\]
whenever $\|u\|_X<\delta$. If $\|h\|_X<\delta$, then $\|th\|_X<\delta$ for every $t\in[0,1]$, and the first estimate gives the final inequality.
\end{proof}

\subsection{Exercises}
\label{subsec:ejercicios-derivadas-superiores}

\begin{exercise}
\label{ejer:derivadas-potencia-multilineal}
Let $A\in\mathcal L_s^m(X;Y)$ and define $P\colon X\longrightarrow Y$ by $P(x):=A[x]^m$. Prove that $P$ is smooth and that, for $0\leq j\leq m$,
\[
D^jP(x)[h_1,\ldots,h_j]
=
\frac{m!}{(m-j)!}A(\underbrace{x,\ldots,x}_{m-j\text{ times}},h_1,\ldots,h_j).
\]
Show that $D^jP=0$ for $j>m$.
\end{exercise}

\begin{exercise}
\label{ejer:derivadas-norma-hilbert}
Let $H$ be a real Hilbert space. Compute the first and second Fréchet derivatives of $u\longmapsto\displaystyle\frac{1}{2}\|u\|_H^2$. For $p>1$, compute the first derivative of $u\longmapsto\displaystyle\frac{1}{p}\|u\|_H^p$ and determine at which points the second derivative exists when $1<p<2$ and when $p\geq2$.
\end{exercise}

\begin{exercise}
\label{ejer:derivadas-composicion-funcional-lineal}
Let $\ell\in X'$ and let $\phi\colon \mathbb R\longrightarrow\mathbb R$ be of class $C^m$. Prove that $f:=\phi\circ\ell$ belongs to $C^m(X,\mathbb R)$ and that
\[
D^jf(x)[h_1,\ldots,h_j]
=
\phi^{(j)}(\ell(x))\prod_{i=1}^{j}\ell(h_i)
\]
for every $1\leq j\leq m$.
\end{exercise}

\begin{exercise}
\label{ejer:energia-dirichlet-haz-derivadas}
Let $(M,\mathbf{g})$ be a compact Riemannian manifold with smooth boundary, and let $\mathbf{E}\to M$ be a smooth vector bundle with bundle metric $\mathbf{h}_{\mathbf{E}}$ and compatible connection $\nabla^{\mathbf{E}}$. Consider
$\displaystyle J(u):=\displaystyle\frac{1}{2}\int_M|\nabla_w^{\mathbf{E}}u|_{\mathbf{g},\mathbf{h}_{\mathbf{E}}}^2\,d\lambda_{\mathbf{g}}$ on $W_0^{1,2}(M,\mathbf{E})$. Prove that $J$ is smooth, compute $DJ(u)$ and $D^2J(u)$, and identify the continuous bilinear operator associated with the second derivative.
\end{exercise}

\begin{exercise}
\label{ejer:taylor-exponencial-banach}
Let $X$ be a unital Banach algebra. Use Taylor's formula to study the exponential map $\displaystyle\exp(x):=\displaystyle\sum_{j=0}^{\infty}\displaystyle\frac{x^j}{j!}$ near the origin. Compute its first derivatives at $0$ and explain the simplification that occurs when the elements under consideration commute.
\end{exercise}

\section{Fundamental theorems of differential calculus}
\label{sec:teoremas-fuertes-calculo-banach}

The following local results turn information about the derivative into information about the function itself. The Banach fixed point theorem~\ref{teo:principio-contraccion-banach} provides the mechanism for existence and uniqueness; the inverse function theorem (Theorem~\ref{teo:funcion-inversa-local-banach}) constructs nonlinear coordinates, and the implicit function theorem (Theorem~\ref{teo:funcion-implicita-banach}) allows some of the variables to be solved for locally. These tools will underlie the adapted charts and regular level sets used on Banach manifolds.

\begin{proposition}[Inversion of operators]
\label{prop:inversion-operadores-suave}
\index{linear operator!inversion}
Let $X$ and $Y$ be Banach spaces, and let $\operatorname{Iso}(X,Y)\subseteq\mathcal L(X,Y)$ be the set of continuous linear isomorphisms. Then $\operatorname{Iso}(X,Y)$ is open, and the map
\[
\operatorname{Inv}\colon \operatorname{Iso}(X,Y)\longrightarrow\operatorname{Iso}(Y,X),
\qquad
\operatorname{Inv}(A):=A^{-1},
\]
is smooth. Its derivative is given by
\begin{equation}
\label{eq:derivada-inversion-operadores}
D\operatorname{Inv}(A)(H)=-A^{-1}HA^{-1}.
\end{equation}
\end{proposition}

\begin{proof}
Let $A\in\operatorname{Iso}(X,Y)$. If $X=\{0\}$, existence of $A$ forces $Y=\{0\}$, and all assertions can be checked directly in the one-point spaces. Suppose $X\neq\{0\}$; then $\|A^{-1}\|_{\mathcal L(Y,X)}>0$. If $H\in\mathcal L(X,Y)$ satisfies $\|A^{-1}H\|_{\mathcal L(X)}<1$, the Neumann series of Theorem~\ref{teo:serie-de-neumann} converges in $\mathcal L(X)$ and
\[
(A+H)^{-1}
=
\displaystyle\sum_{j=0}^{\infty}(-A^{-1}H)^jA^{-1}.
\]
Thus $A+H$ is invertible, and $\operatorname{Iso}(X,Y)$ is open. If $\|A^{-1}H\|_{\mathcal L(X)}\leq\displaystyle\frac{1}{2}$, then
\[
\|(A+H)^{-1}\|_{\mathcal L(Y,X)}
\leq
\frac{\|A^{-1}\|_{\mathcal L(Y,X)}}{1-\|A^{-1}H\|_{\mathcal L(X)}}
\leq
2\|A^{-1}\|_{\mathcal L(Y,X)}.
\]
The identity
\[
(A+H)^{-1}-A^{-1}=-(A+H)^{-1}HA^{-1}
\]
implies
\[
\|(A+H)^{-1}-A^{-1}\|_{\mathcal L(Y,X)}
\leq
2\|A^{-1}\|_{\mathcal L(Y,X)}^2\|H\|_{\mathcal L(X,Y)},
\]
so the inversion map is continuous at $A$.

Moreover,
\[
(A+H)^{-1}-A^{-1}+A^{-1}HA^{-1}
=
A^{-1}HA^{-1}H(A+H)^{-1}.
\]
When $0<\|H\|_{\mathcal L(X,Y)}\leq\displaystyle\frac{1}{2\|A^{-1}\|_{\mathcal L(Y,X)}}$, the preceding estimate gives
\[
\frac{\|(A+H)^{-1}-A^{-1}+A^{-1}HA^{-1}\|_{\mathcal L(Y,X)}}
{\|H\|_{\mathcal L(X,Y)}}
\leq
2\|A^{-1}\|_{\mathcal L(Y,X)}^3\|H\|_{\mathcal L(X,Y)}.
\]
The right-hand side tends to zero as $H\to0$, proving \eqref{eq:derivada-inversion-operadores}.

The derivative formula can be written as
$D\operatorname{Inv}(A)(H)=-\operatorname{Inv}(A)H\operatorname{Inv}(A)$. The map assigning to $(B,C)\in\mathcal L(Y,X)\times\mathcal L(Y,X)$ the operator $H\longmapsto-BHC$ is continuous and bilinear with values in $\mathcal L(\mathcal L(X,Y),\mathcal L(Y,X))$, since its norm is at most $\|B\|\|C\|$. By Proposition~\ref{prop:diferenciabilidad-operadores-multilineales}, this map is smooth. Since $\operatorname{Inv}$ is continuous, its derivative formula first proves continuity of $D\operatorname{Inv}$; hence $\operatorname{Inv}$ is of class $C^1$. Now suppose it is of class $C^r$ for some integer $r\geq1$. Composing the preceding bilinear map with $A\longmapsto(\operatorname{Inv}(A),\operatorname{Inv}(A))$, Corollary~\ref{cor:composicion-Ck-banach} shows that $D\operatorname{Inv}$ is of class $C^r$. Consequently, $\operatorname{Inv}$ is of class $C^{r+1}$. This induction step proves smoothness.
\end{proof}

\begin{definition}[Diffeomorphism between open subsets of Banach spaces]
\label{def:difeomorfismo-banach}
\index{diffeomorphism!between open subsets of Banach spaces}
Let $U\subseteq X$ and $V\subseteq Y$ be open subsets of Banach spaces. A bijection $f\colon U\longrightarrow V$ is a \textbf{diffeomorphism of class $C^k$} if $f\in C^k(U,Y)$ and $f^{-1}\in C^k(V,X)$.
\end{definition}

\begin{theorem}[Local inverse function theorem]
\label{teo:funcion-inversa-local-banach}
\index{inverse function theorem!in Banach spaces}
Let $X$ and $Y$ be Banach spaces, let $U\subseteq X$ be open, and let $f\in C^1(U,Y)$. If $a\in U$ and $Df(a)\colon X\longrightarrow Y$ is an isomorphism, there exist open neighborhoods $U_0\subseteq U$ of $a$ and $V_0\subseteq Y$ of $f(a)$ such that $f\restriction_{U_0}\colon U_0\longrightarrow V_0$ is a diffeomorphism of class $C^1$. If $g\colon V_0\longrightarrow U_0$ is its inverse, then
\begin{equation}
\label{eq:derivada-inversa-banach}
Dg(y)=Df(g(y))^{-1}
\end{equation}
for every $y\in V_0$. If $f\in C^k(U,Y)$, where $k\in\mathbb N$, $k\geq1$, or $k=\infty$, then $g\in C^k(V_0,X)$.
\end{theorem}

\begin{proof}
If $X=\{0\}$, invertibility of $Df(a)$ implies $Y=\{0\}$, and the result is immediate. Suppose $X\neq\{0\}$. Write $A:=Df(a)^{-1}\in\mathcal L(Y,X)$. By continuity of $Df$ at $a$, there exists $r>0$ such that $\overline{B}_{d_X}(a,2r)\subseteq U$ and
\begin{equation}
\label{eq:control-derivada-inversa-local}
\|A(Df(x)-Df(a))\|_{\mathcal L(X)}\leq\frac{1}{2}
\end{equation}
for every $x\in B_{d_X}(a,2r)$. For each $y\in Y$, define
\[
T_y\colon \overline{B}_{d_X}(a,r)\longrightarrow X,
\qquad
T_y(x):=x-A(f(x)-y).
\]
Let $x_1,x_2\in\overline{B}_{d_X}(a,r)$. The segment $[x_1,x_2]$ is contained in $B_{d_X}(a,2r)$. The map $H(x):=x-Af(x)$ is of class $C^1$ and, by \eqref{eq:control-derivada-inversa-local},
\[
\|DH(x)\|_{\mathcal L(X)}
=
\|I_X-ADf(x)\|_{\mathcal L(X)}
=
\|A(Df(a)-Df(x))\|_{\mathcal L(X)}
\leq
\frac{1}{2}.
\]
Theorem~\ref{teo:desigualdades-valor-medio-direccionales} implies
\begin{equation}
\label{eq:contraccion-funcion-inversa}
\|T_y(x_1)-T_y(x_2)\|_X
\leq
\frac{1}{2}\|x_1-x_2\|_X.
\end{equation}

Let $b:=f(a)$ and choose $\delta>0$ such that $\|A\|_{\mathcal L(Y,X)}\delta<\displaystyle\frac{r}{2}$. If $y\in B_{d_Y}(b,\delta)$, then
\[
\|T_y(a)-a\|_X
=
\|A(y-b)\|_X
<
\frac{r}{2}.
\]
For each $x\in\overline{B}_{d_X}(a,r)$, inequality \eqref{eq:contraccion-funcion-inversa} gives
\[
\|T_y(x)-a\|_X
\leq
\|T_y(x)-T_y(a)\|_X+\|T_y(a)-a\|_X
<
\frac{1}{2}\|x-a\|_X+\frac{r}{2}
\leq r.
\]
Thus $T_y$ maps $\overline{B}_{d_X}(a,r)$ into itself and is a contraction. Theorem~\ref{teo:principio-contraccion-banach} provides a unique fixed point $g(y)\in\overline{B}_{d_X}(a,r)$. The equality $T_y(g(y))=g(y)$ is equivalent to $f(g(y))=y$. Moreover,
\[
\|g(y)-a\|_X
\leq
\frac{1}{1-\frac{1}{2}}\|T_y(a)-a\|_X
<r,
\]
so $g(y)\in B_{d_X}(a,r)$.

If $y_1,y_2\in B_{d_Y}(b,\delta)$, then, using \eqref{eq:contraccion-funcion-inversa},
\[
\begin{aligned}
\|g(y_1)-g(y_2)\|_X
&=\|T_{y_1}(g(y_1))-T_{y_2}(g(y_2))\|_X\\
&\leq\|T_{y_1}(g(y_1))-T_{y_1}(g(y_2))\|_X
+\|A(y_1-y_2)\|_X\\
&\leq\frac{1}{2}\|g(y_1)-g(y_2)\|_X
+\|A\|_{\mathcal L(Y,X)}\|y_1-y_2\|_Y.
\end{aligned}
\]
Consequently,
\begin{equation}
\label{eq:lipschitz-inversa-local}
\|g(y_1)-g(y_2)\|_X
\leq
2\|A\|_{\mathcal L(Y,X)}\|y_1-y_2\|_Y.
\end{equation}
In particular, $g$ is continuous.

Define $V_0:=B_{d_Y}(b,\delta)$ and $U_0:=f^{-1}(V_0)\cap B_{d_X}(a,r)$. If $x\in U_0$ and $y=f(x)$, then $x$ is a fixed point of $T_y$ in $\overline{B}_{d_X}(a,r)$; by uniqueness of the fixed point, $x=g(y)$. Thus $f\restriction_{U_0}$ is bijective onto $V_0$, with inverse $g$.

For each $x\in U_0$, the identity
\[
ADf(x)=I_X+A(Df(x)-Df(a))
\]
and \eqref{eq:control-derivada-inversa-local} show, through Theorem~\ref{teo:serie-de-neumann}, that $ADf(x)$ is invertible. Thus $Df(x)$ is an isomorphism.

Fix $y\in V_0$ and write $x:=g(y)$. Let $k\in Y$ be such that $y+k\in V_0$ and define $h:=g(y+k)-g(y)$. Then $k=f(x+h)-f(x)$ and
\begin{equation}
\label{eq:resto-derivada-inversa-local}
g(y+k)-g(y)-Df(x)^{-1}k
=
-Df(x)^{-1}\bigl(f(x+h)-f(x)-Df(x)h\bigr).
\end{equation}
Let $\varepsilon>0$. By Fréchet differentiability of $f$ at $x$, there exists $\eta>0$ such that
\[
\|f(x+u)-f(x)-Df(x)u\|_Y
<
\frac{\varepsilon}{2\bigl(1+\|A\|_{\mathcal L(Y,X)}\bigr)\bigl(1+\|Df(x)^{-1}\|_{\mathcal L(Y,X)}\bigr)}\|u\|_X
\]
whenever $0<\|u\|_X<\eta$. By \eqref{eq:lipschitz-inversa-local}, we have $\|h\|_X\leq2\|A\|_{\mathcal L(Y,X)}\|k\|_Y$. Consequently, if
\[
0<\|k\|_Y<\frac{\eta}{2\|A\|_{\mathcal L(Y,X)}},
\]
then $\|h\|_X<\eta$, so \eqref{eq:resto-derivada-inversa-local} implies
\[
\frac{\|g(y+k)-g(y)-Df(x)^{-1}k\|_X}{\|k\|_Y}
<\varepsilon.
\]
This proves that $g$ is Fréchet differentiable at $y$ and that $Dg(y)=Df(g(y))^{-1}$.

Continuity of $Dg$ follows from continuity of $g$, continuity of $Df$, and continuity of the inversion map established in Proposition~\ref{prop:inversion-operadores-suave}. Thus $g\in C^1(V_0,X)$. If $f\in C^k(U,Y)$ and $k>1$, start with the $C^1$ regularity of $g$ just obtained. Suppose $g$ is of class $C^j$, where $j\in\mathbb N$ and $1\leq j<k$. Since $Df$ is of class $C^{k-1}$, it is in particular of class $C^j$. Corollary~\ref{cor:composicion-Ck-banach} and smoothness of the inversion map imply that $\operatorname{Inv}\circ Df\circ g$ is of class $C^j$. By \eqref{eq:derivada-inversa-banach}, this is precisely $Dg$; hence $g$ is of class $C^{j+1}$. Induction reaches $j=k$ when $k$ is finite. If $k=\infty$, the same argument applies at every finite order and proves that $g$ is smooth.
\end{proof}

Invertibility of the derivative is a local condition and alone does not guarantee global injectivity. The uniform assumption in the following result prevents local inverses from degenerating as they are continued along paths.

\begin{lemma}[Path lifting]
\label{lem:levantamiento-caminos-inversa-global}
Let $X$ and $Y$ be Banach spaces, and let $f\in C^1(X,Y)$. Suppose $Df(x)$ is an isomorphism for every $x\in X$ and there exists $C>0$ such that $\|Df(x)^{-1}\|_{\mathcal L(Y,X)}\leq C$ for every $x\in X$. Let $\gamma\colon [0,1]\longrightarrow Y$ be of class $C^1$, and let $x_0\in X$ be such that $f(x_0)=\gamma(0)$. Then there exists a unique map $\widetilde\gamma\colon [0,1]\longrightarrow X$ of class $C^1$ such that $\widetilde\gamma(0)=x_0$ and $f\circ\widetilde\gamma=\gamma$. Moreover,
\begin{equation}
\label{eq:cota-derivada-levantamiento}
\|\widetilde\gamma'(t)\|_X\leq C\|\gamma'(t)\|_Y
\end{equation}
for every $t\in[0,1]$.
\end{lemma}

\begin{proof}
Theorem~\ref{teo:funcion-inversa-local-banach}, applied at $x_0$, gives a unique lift on an interval $[0,\tau_0]$ with $\tau_0>0$. Let
\[
\mathcal A:=\{\tau\in[0,1]\mid \gamma\restriction_{[0,\tau]}\text{ admits a lift of class }C^1\text{ starting at }x_0\}
\]
and let $\displaystyle \alpha:=\displaystyle\sup\mathcal A$. Let us establish uniqueness on common intervals precisely. If two continuous lifts have the same initial value, the set of times when they agree is nonempty and closed by continuity. It is also open in the common interval: at a time of agreement, both paths remain, for a small interval, in a common neighborhood on which $f$ is injective. Since they have the same image under $f$, they agree there. Connectedness of the interval implies that they agree throughout it. Thus the partial lifts glue together to define a unique map $\widetilde\gamma\colon [0,\alpha)\longrightarrow X$, of class $C^1$ on every compact interval contained in its domain.

Differentiating $f(\widetilde\gamma(t))=\gamma(t)$ and applying the chain rule gives
\[
Df(\widetilde\gamma(t))\widetilde\gamma'(t)=\gamma'(t),
\qquad
\widetilde\gamma'(t)=Df(\widetilde\gamma(t))^{-1}\gamma'(t),
\]
which proves \eqref{eq:cota-derivada-levantamiento}. If $0\leq t_1<t_2<\alpha$, the integral mean value formula gives
\[
\|\widetilde\gamma(t_2)-\widetilde\gamma(t_1)\|_X
\leq
C\left(\max_{t\in[0,1]}\|\gamma'(t)\|_Y\right)|t_2-t_1|.
\]
This estimate shows that $\widetilde\gamma(t)$ is Cauchy as $t\to\alpha^{-}$. By completeness of $X$, there exists $\displaystyle x_\alpha:=\lim_{t\to\alpha^{-}}\widetilde\gamma(t)$, whether $\alpha<1$ or $\alpha=1$. Define $\widetilde\gamma(\alpha):=x_\alpha$. The extension is continuous and satisfies $f(x_\alpha)=\gamma(\alpha)$.

Choose a local inverse $g_\alpha\colon V_\alpha\longrightarrow U_\alpha$ of $f$, with $x_\alpha\in U_\alpha$ and $\gamma(\alpha)\in V_\alpha$. By the convergence just proved and continuity of $\gamma$, there exists $\varepsilon>0$ such that $\widetilde\gamma(t)\in U_\alpha$ and $\gamma(t)\in V_\alpha$ for $t\in(\alpha-\varepsilon,\alpha]$. On this interval, injectivity of $f\restriction_{U_\alpha}$ gives
\[
\widetilde\gamma(t)=g_\alpha(\gamma(t)).
\]
Thus the extension is of class $C^1$ up to $\alpha$, with a one-sided derivative at the endpoint. If $\alpha<1$, continuity of $\gamma$ allows this same formula to be used on $[\alpha,\alpha+\varepsilon']$ for some $\varepsilon'>0$ with $\alpha+\varepsilon'\leq1$. We would obtain a lift of class $C^1$ beyond $\alpha$, contradicting the definition of the supremum. Hence $\alpha=1$.

In particular, the value $\widetilde\gamma(1)$ is already defined, and the formula involving $g_\alpha$ proves $C^1$ regularity at the final endpoint. Differentiating it gives $\widetilde\gamma'(1)=Df(\widetilde\gamma(1))^{-1}\gamma'(1)$ there as well, so \eqref{eq:cota-derivada-levantamiento} holds on all of $[0,1]$. Uniqueness on the closed interval is the uniqueness established at the beginning of the proof.
\end{proof}

\begin{lemma}[Homotopy lifting]
\label{lem:levantamiento-homotopias-inversa-global}
Under the assumptions of the preceding lemma, let $H\colon [0,1]^2\longrightarrow Y$ be of class $C^1$, and let $\eta\colon [0,1]\longrightarrow X$ be a continuous path such that $f(\eta(s))=H(0,s)$ for every $s\in[0,1]$. Then there exists a unique continuous map $\widetilde H\colon [0,1]^2\longrightarrow X$ such that $f\circ\widetilde H=H$ and $\widetilde H(0,s)=\eta(s)$. For each $s$, the path $t\longmapsto\widetilde H(t,s)$ is of class $C^1$.
\end{lemma}

\begin{proof}
For each $s\in[0,1]$, Lemma~\ref{lem:levantamiento-caminos-inversa-global} gives a unique lift of the path $t\longmapsto H(t,s)$ starting at $\eta(s)$; denote it by $t\longmapsto\widetilde H(t,s)$. We must prove that the resulting map is continuous in both variables.

Let
\[
\mathcal T:=\{\tau\in[0,1]\mid \widetilde H\text{ is continuous on }[0,\tau]\times[0,1]\}.
\]
We will prove that $\mathcal T=[0,1]$. For each $s_0\in[0,1]$, Theorem~\ref{teo:funcion-inversa-local-banach} gives open neighborhoods $U_{s_0}$ of $\eta(s_0)$ and $V_{s_0}$ of $H(0,s_0)$ such that $f\restriction_{U_{s_0}}\colon U_{s_0}\longrightarrow V_{s_0}$ is a homeomorphism. By continuity of $\eta$ and $H$, there exist an open interval $I_{s_0}\subseteq[0,1]$ containing $s_0$ and a number $\delta_{s_0}>0$ such that $\eta(s)\in U_{s_0}$ and $H(t,s)\in V_{s_0}$ for $(t,s)\in[0,\delta_{s_0}]\times I_{s_0}$. On this rectangle, uniqueness of lifts implies
\[
\widetilde H(t,s)=(f\restriction_{U_{s_0}})^{-1}(H(t,s)).
\]
Compactness of $[0,1]$ allows us to choose finitely many intervals $I_{s_0}$ covering it. If $\delta$ is the minimum of the corresponding numbers $\delta_{s_0}$, the preceding formula proves that $\widetilde H$ is continuous on $[0,\delta]\times[0,1]$. Thus $\mathcal T$ is nonempty.

Let $\tau\in\mathcal T$ with $\tau<1$. The set $\widetilde H(\{\tau\}\times[0,1])$ is compact. For each $s_0\in[0,1]$, apply the local theorem again at $\widetilde H(\tau,s_0)$. Continuity of $s\longmapsto\widetilde H(\tau,s)$ and $H$ allows us to find an interval $I_{s_0}$ and a number $\delta_{s_0}>0$ such that, for $s\in I_{s_0}$ and $t\in[\tau,\tau+\delta_{s_0}]\cap[0,1]$, the lift is given by a common local inverse. A finite subcover of $[0,1]$ and the minimum of the numbers $\delta_{s_0}$ show that there exists $\delta>0$ such that $\widetilde H$ is continuous on $\displaystyle [0,\min\{1,\tau+\delta\}]\times[0,1]$. Thus $\mathcal T$ is open to the right.

To check that continuation cannot stop, let $\displaystyle \tau_*:=\displaystyle\sup\mathcal T$. Since $H$ is of class $C^1$ on the compact square, there exists
\[
M:=\max_{(t,s)\in[0,1]^2}\left\|\frac{\partial H}{\partial t}(t,s)\right\|_Y<\infty.
\]
For each $s\in[0,1]$ and $0\leq t_1<t_2<\tau_*$, estimate \eqref{eq:cota-derivada-levantamiento} and the integral mean value formula give
\[
\|\widetilde H(t_2,s)-\widetilde H(t_1,s)\|_X
\leq
CM|t_2-t_1|.
\]
The estimate is uniform in $s$. Consequently, as $t\to \tau_*^{-}$, $\widetilde H(t,s)$ has a uniform limit in $s$. For each fixed $s$, the individual lift is already continuous on $[0,1]$, so this limit agrees with the value $\widetilde H(\tau_*,s)$ defined at the beginning of the proof. Moreover, for every $t<\tau_*$, there exists an element of $\mathcal T$ greater than $t$, so $s\longmapsto\widetilde H(t,s)$ is continuous. The uniform limit of these functions is continuous. We thus obtain continuity of $s\longmapsto\widetilde H(\tau_*,s)$ and, by passage to the limit in the estimate,
\[
\|\widetilde H(t,s)-\widetilde H(\tau_*,s_0)\|_X
\leq CM(\tau_*-t)+\|\widetilde H(\tau_*,s)-\widetilde H(\tau_*,s_0)\|_X
\]
for $t\in[0,\tau_*]$ and $s,s_0\in[0,1]$. The right-hand side tends to zero as $(t,s)\longrightarrow(\tau_*,s_0)$ with $t\leq\tau_*$. This proves joint continuity up to this edge, that is, $\tau_*\in\mathcal T$. The local argument in the preceding paragraph would extend continuity beyond $\tau_*$ if $\tau_*<1$, contradicting the definition of the supremum. Thus $\tau_*=1$, and $\widetilde H$ is continuous on the entire square.

The identity $f\circ\widetilde H=H$ and the initial condition hold by construction. Uniqueness follows for each fixed value of $s$ from the argument for continuous lifts in the proof of Lemma~\ref{lem:levantamiento-caminos-inversa-global}.
\end{proof}

\begin{theorem}[Global inverse function theorem]
\label{teo:funcion-inversa-global-banach}
\index{inverse function theorem!global}
Let $X$ and $Y$ be Banach spaces, and let $f\in C^1(X,Y)$. Suppose $Df(x)$ is an isomorphism for every $x\in X$ and there exists $C>0$ such that $\|Df(x)^{-1}\|_{\mathcal L(Y,X)}\leq C$ for every $x\in X$. Then $f\colon X\longrightarrow Y$ is a global diffeomorphism of class $C^1$. If, in addition, $f\in C^k(X,Y)$, where $k\in\mathbb N$, $k\geq1$, or $k=\infty$, its inverse is of class $C^k$.
\end{theorem}

\begin{proof}
Fix $a\in X$ and let $b:=f(a)$. Given $y\in Y$, consider the segment $\gamma(t):=(1-t)b+ty$. By Lemma~\ref{lem:levantamiento-caminos-inversa-global}, there exists a lift $\widetilde\gamma\colon [0,1]\longrightarrow X$ with $\widetilde\gamma(0)=a$. Then $f(\widetilde\gamma(1))=\gamma(1)=y$. Since $y$ was arbitrary, $f$ is surjective.

To prove injectivity, suppose $f(x_0)=f(x_1)=y_0$. Let $\sigma(t):=(1-t)x_0+tx_1$ and let $\gamma:=f\circ\sigma$. The map
\[
H(t,s):=(1-s)\gamma(t)+sy_0
\]
is a homotopy between the loop $\gamma$ and the constant loop $y_0$ and satisfies $H(0,s)=y_0$ for every $s$. Apply Lemma~\ref{lem:levantamiento-homotopias-inversa-global} with $\eta(s):=x_0$. For each $s\in[0,1]$, the point $\widetilde H(1,s)$ belongs to $f^{-1}(y_0)$. This fiber is discrete because $f$ is locally injective at each of its points. Since $s\longmapsto\widetilde H(1,s)$ is continuous and $[0,1]$ is connected, this map is constant. For $s=0$, uniqueness of the lift implies $\widetilde H(t,0)=\sigma(t)$ and hence $\widetilde H(1,0)=x_1$. For $s=1$, the lift is the constant path $x_0$, so $\widetilde H(1,1)=x_0$. We conclude that $x_0=x_1$.

Therefore, $f$ is bijective. On each neighborhood where
Theorem~\ref{teo:funcion-inversa-local-banach} provides a local inverse,
that inverse agrees with $f^{-1}$ by uniqueness of the preimage. Thus
$f^{-1}$ is of class $C^1$ and satisfies
$D(f^{-1})(y)=Df(f^{-1}(y))^{-1}$. If $f\in C^k$, the same local inverses
are of class $C^k$; since this regularity is local, $f^{-1}\in C^k$
on all of $Y$. For $k=\infty$, the assertion holds at every finite order.
\end{proof}

\begin{theorem}[Implicit function theorem]
\label{teo:funcion-implicita-banach}
\index{implicit function theorem!in Banach spaces}
Let $X$, $Y$, and $Z$ be Banach spaces, let $U\subseteq X\times Y$ be open, and let $F\in C^k(U,Z)$, where $k\in\mathbb N$, $k\geq1$, or $k=\infty$. Suppose $(a,b)\in U$, $F(a,b)=0$, and $D_2F(a,b)\colon Y\longrightarrow Z$ is an isomorphism. Then there exist open neighborhoods $U_0\subseteq X$ of $a$ and $V_0\subseteq Y$ of $b$, and a unique function $\varphi\colon U_0\longrightarrow V_0$ of class $C^k$ such that, for every $(x,y)\in U_0\times V_0$,
\begin{equation}
\label{eq:caracterizacion-funcion-implicita}
F(x,y)=0\quad\Longleftrightarrow\quad y=\varphi(x).
\end{equation}
After shrinking $U_0$ and $V_0$ if necessary, $D_2F(x,\varphi(x))$ is an isomorphism for every $x\in U_0$ and
\begin{equation}
\label{eq:derivada-funcion-implicita}
D\varphi(x)
=
-D_2F(x,\varphi(x))^{-1}\circ D_1F(x,\varphi(x)).
\end{equation}
\end{theorem}

\begin{proof}
Define
\[
\Phi\colon U\longrightarrow X\times Z,
\qquad
\Phi(x,y):=\bigl(x,F(x,y)\bigr).
\]
Proposition~\ref{prop:derivada-total-implica-parciales} implies that
\[
D\Phi(a,b)(h,k)
=
\bigl(h,D_1F(a,b)h+D_2F(a,b)k\bigr).
\]
This operator is an isomorphism from $X\times Y$ onto $X\times Z$. Its inverse is given by
\[
(u,v)
\longmapsto
\left(u,D_2F(a,b)^{-1}\bigl(v-D_1F(a,b)u\bigr)\right).
\]
Theorem~\ref{teo:funcion-inversa-local-banach} gives open neighborhoods $\mathcal U\subseteq U$ of $(a,b)$ and $\mathcal V\subseteq X\times Z$ of $(a,0)$ such that $\Phi\restriction_{\mathcal U}\colon \mathcal U\longrightarrow\mathcal V$ is a diffeomorphism of class $C^k$. First choose an open neighborhood $V_0\subseteq Y$ of $b$ such that there exists an open neighborhood $U_1\subseteq X$ of $a$ with $U_1\times V_0\subseteq\mathcal U$. By continuity of $\Phi^{-1}$ at $(a,0)$, there exist open neighborhoods $U_0\subseteq U_1$ of $a$ and $W_0\subseteq Z$ of $0$ such that $U_0\times W_0\subseteq\mathcal V$ and
\[
\Phi^{-1}(U_0\times W_0)\subseteq U_1\times V_0.
\]

For each $x\in U_0$, define $(x,\varphi(x)):=\Phi^{-1}(x,0)$. The first component of $\Phi$ is the identity on $X$, so the first component of $\Phi^{-1}(x,0)$ must be $x$, and the preceding choice ensures that $\varphi(x)\in V_0$. Moreover,
\[
\Phi(x,\varphi(x))=(x,0)
\quad\Longleftrightarrow\quad
F(x,\varphi(x))=0.
\]
If $(x,y)\in U_0\times V_0$ and $F(x,y)=0$, then $\Phi(x,y)=(x,0)$; injectivity of $\Phi\restriction_{\mathcal U}$ implies $y=\varphi(x)$. This proves \eqref{eq:caracterizacion-funcion-implicita}. Since $\varphi$ is obtained by composing the map $x\longmapsto(x,0)$, the inverse $\Phi^{-1}$, and the projection onto $Y$, it belongs to $C^k(U_0,Y)$.

Continuity of $D_2F$ and openness of $\operatorname{Iso}(Y,Z)$, proved in Proposition~\ref{prop:inversion-operadores-suave}, allow the neighborhoods to be shrunk so that $D_2F(x,\varphi(x))$ is invertible for every $x\in U_0$. Differentiating the identity $F(x,\varphi(x))=0$ and applying the chain rule and \eqref{eq:descomposicion-derivada-total-frechet} gives
\[
D_1F(x,\varphi(x))+D_2F(x,\varphi(x))\circ D\varphi(x)=0.
\]
Composition on the left with $D_2F(x,\varphi(x))^{-1}$ yields \eqref{eq:derivada-funcion-implicita}.
\end{proof}

\section{Examples of calculus in function spaces}
\label{sec:ejemplos-calculo-algebras-funciones}

Pointwise operations on functions provide examples of maps between Banach
spaces. To apply the preceding calculus, their regularity must be verified
in the norm of the space, in addition to their regularity at each point.
We first establish the estimates that will be used to construct manifolds
of functions and curves.

By a real commutative unital Banach algebra we mean a Banach space
$\mathcal A$ with a bilinear, associative, and commutative product,
an element $1$ satisfying $1u=u1=u$, and a norm satisfying
$\|uv\|_{\mathcal A}\leq\|u\|_{\mathcal A}\|v\|_{\mathcal A}$.
The unit may have norm different from one.

\begin{lemma}[The algebra of functions of finite regularity]
\label{lem:algebra-ck-compacta}
Let $(M,\mathbf g)$ be a smooth, nonempty compact Riemannian manifold without
boundary, and let $k\in\mathbb N_0$. The space $C^k(M,\mathbb R)$ is a
Banach space with the norm
\[
\|f\|_{k,\infty}:=\sum_{j=0}^{k}\sup_{x\in M}|\nabla^jf(x)|_{\mathbf g}.
\]
The equivalent norm $\|f\|_{\mathcal A}:=2^k\|f\|_{k,\infty}$ satisfies
\[
\|fg\|_{\mathcal A}\leq\|f\|_{\mathcal A}\|g\|_{\mathcal A},
\qquad
\|f\|_{C^0(M)}\leq\|f\|_{\mathcal A}.
\]
With pointwise operations and the constant function $1$ as its unit,
$\mathcal A:=C^k(M,\mathbb R)$ is therefore a commutative Banach algebra.
The spaces $C^k(M,\mathbf E)$ of sections of a smooth vector bundle of
finite rank, with a smooth metric and connection, are also complete with
the analogous norm. If this norm is multiplied by $2^k$, multiplication
of a section by a function satisfies
$\|f\mathbf s\|_{\mathcal A,\mathbf E}
\leq\|f\|_{\mathcal A}\|\mathbf s\|_{\mathcal A,\mathbf E}$.
\end{lemma}

\begin{proof}
Let $(f_\nu)$ be a Cauchy sequence for $\|\cdot\|_{k,\infty}$.
For each $j\in\{0,\ldots,k\}$, the continuous sections
$\nabla^jf_\nu$ converge uniformly to a continuous section $F_j$:
in a trivialization, their components are uniformly Cauchy sequences of
continuous functions, and their limits agree on overlaps by the
transition formulas. In particular, $f_\nu$ converges uniformly to
$f:=F_0$.

We verify that $F_j=\nabla^jf$. In a coordinate ball whose closure is
compact and contained in a chart, the components of $\nabla^jf_\nu$
are the partial derivatives of order $j$ of $f_\nu$ plus linear
combinations of derivatives of lower order, with fixed smooth
coefficients. Conversely, solving in increasing order, each partial
derivative of order $j\leq k$ is a combination of components of
$f_\nu,\nabla f_\nu,\ldots,\nabla^jf_\nu$, with fixed smooth
coefficients. These coefficients are bounded on the closed ball.
Consequently, all partial derivatives of order at most $k$ converge
uniformly there.

If $G_\alpha$ is the limit of $\partial^\alpha f_\nu$, the fundamental
theorem of calculus along a coordinate segment contained in the ball gives
\[
\partial^\alpha f_\nu(x+te_i)-\partial^\alpha f_\nu(x)
=\int_0^t\partial^{\alpha+e_i}f_\nu(x+se_i)\,ds,
\qquad |\alpha|<k.
\]
The absolute value of the error upon replacing the integrand by its
limit is at most
$|t|\displaystyle\sup_s|\partial^{\alpha+e_i}f_\nu(x+se_i)-G_{\alpha+e_i}(x+se_i)|$,
which tends to zero. Thus
\[
G_\alpha(x+te_i)-G_\alpha(x)
=\int_0^tG_{\alpha+e_i}(x+se_i)\,ds.
\]
Continuity of $G_{\alpha+e_i}$ shows that
$\partial_iG_\alpha=G_{\alpha+e_i}$. Induction on the order of the
derivatives yields $f\in C^k$ and $F_j=\nabla^jf$. The coordinate
balls cover $M$, so the identification is global and
$\|f_\nu-f\|_{k,\infty}\to0$. For $k=0$, the initial uniform
convergence suffices. The same argument, applied to the components of
a section and the local formulas for its connection, proves completeness
of $C^k(M,\mathbf E)$.

After applying the Leibniz rule for the connection $j$ times, each term
contains a derivative of order $a$ of the first factor and one of order
$j-a$ of the second. There are $\binom ja$ terms of this type;
permutations of the covariant arguments preserve the tensor norm.
Therefore,
\[
\begin{split}
\|fg\|_{k,\infty}
&\leq\sum_{j=0}^{k}\sum_{a=0}^{j}\binom ja
\|\nabla^af\|_\infty\|\nabla^{j-a}g\|_\infty\\
&\leq2^k
\left(\sum_{a=0}^{k}\|\nabla^af\|_\infty\right)
\left(\sum_{b=0}^{k}\|\nabla^bg\|_\infty\right).
\end{split}
\]
Multiplying by $2^k$ gives the stated inequality. The same calculation,
with a vector-valued second factor, proves the estimate for
$f\mathbf s$. The uniform bound on $f$ is one of the terms of the norm.
\end{proof}

\begin{lemma}[Products and evaluation in a Sobolev space of curves]
\label{lem:algebra-h1-intervalo}
Each class in $H^1((0,1),\mathbb R)=W^{1,2}((0,1),\mathbb R)$ has a
unique absolutely continuous representative on $[0,1]$. For this
representative,
\[
\|u\|_{C^0([0,1])}\leq\sqrt2\|u\|_{H^1},
\qquad
\|u\|_{H^1}:=
\left(\int_0^1|u(t)|^2\,dt+\int_0^1|u'(t)|^2\,dt\right)^{\frac{1}{2}}.
\]
The product of two elements of $H^1((0,1),\mathbb R)$ belongs to the
same space and satisfies
\[
(uv)'=u'v+uv',
\qquad
\|uv\|_{H^1}\leq3\sqrt2\|u\|_{H^1}\|v\|_{H^1}.
\]
Consequently, the norm
$\|u\|_{\mathcal A_1}:=3\sqrt2\|u\|_{H^1}$ makes this space a
commutative unital Banach algebra and dominates the uniform norm.
Evaluation at each point of $[0,1]$, including the endpoints, is linear
and continuous. The assertions about representatives and evaluation
apply componentwise to $H^1((0,1),\mathbb R^d)$.
\end{lemma}

\begin{proof}
Let $w\in L^2((0,1))$ be the weak derivative of $u$. Define
$F(t):=\int_0^tw(s)\,ds$. For disjoint intervals $(a_i,b_i)$,
\[
\sum_i|F(b_i)-F(a_i)|
\leq\int_{\bigcup_i(a_i,b_i)}|w(s)|\,ds
\leq\|w\|_{L^2}\left(\sum_i(b_i-a_i)\right)^{\frac{1}{2}}.
\]
This inequality proves absolute continuity of $F$ directly.
Moreover, for $\phi\in C_c^\infty((0,1))$, Fubini gives
\[
\int_0^1F(t)\phi'(t)\,dt
=\int_0^1w(s)\left(\int_s^1\phi'(t)\,dt\right)ds
=-\int_0^1w(s)\phi(s)\,ds.
\]
Fubini applies because the absolute integral of the double integrand
is at most $\|\phi'\|_\infty\|w\|_{L^1}<\infty$.
Thus the weak derivative of $F$ is $w$, and $u-F$ has zero weak
derivative and is constant almost everywhere by
Lemma~\ref{lem:derivadas-debiles-nulas-euclidiano}. This gives the
representative $\widetilde u(t)=c+F(t)$. Two continuous representatives
of the same class agree: if their difference were nonzero at a point,
continuity would keep it nonzero on an interval of positive length.
From now on we use this representative and omit the tilde.

There exists $s\in[0,1]$ such that $|u(s)|\leq\|u\|_{L^2}$, since
the opposite strict inequality throughout the interval contradicts the
integral defining $\|u\|_{L^2}^2$. For every $t\in[0,1]$,
Cauchy--Schwarz gives
\[
|u(t)|\leq|u(s)|+\int_0^1|u'(r)|\,dr
\leq\|u\|_{L^2}+\|u'\|_{L^2}
\leq\sqrt2\|u\|_{H^1}.
\]
To justify the product rule in the Sobolev norm, we construct smooth
approximations. By Proposition~\ref{densidad de Cc en Lp}, choose
$w_\nu\in C_c^\infty((0,1))$ with $w_\nu\to w$ in $L^2$ and set
$u_\nu(t):=c+\int_0^tw_\nu(s)\,ds$. Then $u_\nu$ is smooth up to
the endpoints and
\[
\|u_\nu-u\|_\infty\leq\|w_\nu-w\|_{L^2},\qquad
\|u_\nu-u\|_{H^1}\leq\sqrt2\|w_\nu-w\|_{L^2}\longrightarrow0.
\]
Construct $v_\nu$ in the same way for $v$.
The products $u_\nu v_\nu$ converge uniformly to $uv$.
For their derivatives, for example,
\[
\|u_\nu'v_\nu-u'v\|_{L^2}
\leq\|u_\nu'-u'\|_{L^2}\|v_\nu\|_\infty
+\|u'\|_{L^2}\|v_\nu-v\|_\infty\longrightarrow0,
\]
since $v_\nu$ converges uniformly and its uniform norms are bounded.
The other term satisfies the same inequality with $u$ and $v$
interchanged. Passing to the limit in integration by parts for the
smooth products shows that $uv$ has weak derivative $u'v+uv'$.
Moreover,
\[
\begin{split}
\|uv\|_{L^2}&\leq\sqrt2\|u\|_{H^1}\|v\|_{H^1},\\
\|u'v+uv'\|_{L^2}&\leq
\sqrt2\|u\|_{H^1}\|v\|_{H^1}
+\sqrt2\|u\|_{H^1}\|v\|_{H^1}.
\end{split}
\]
Since the $H^1$ norm is at most the sum of the two $L^2$
norms, the constant $3\sqrt2$ follows. Completeness is that of Sobolev
spaces, already proved. Multiplying the norm by $3\sqrt2$ makes it
submultiplicative, and the preceding bounds prove continuity of
evaluation.
\end{proof}

\begin{lemma}[Power series in a Banach algebra]
\label{lem:series-potencias-algebra-banach}
Let $\mathcal A$ be a real, commutative unital Banach algebra whose norm
is submultiplicative. Suppose that $\displaystyle\sum_{n=0}^\infty a_nt^n$ has
radius of convergence $R>0$. Then
\[
F(u):=\sum_{n=0}^\infty a_nu^n,
\qquad \|u\|_{\mathcal A}<R,
\]
is a smooth map with values in $\mathcal A$, where $u^0:=1$.
Its derivatives are obtained by differentiating the series term by term,
with uniform convergence on every ball of radius $r<R$ in the norm
of the corresponding multilinear operators.
\end{lemma}

\begin{proof}
Set $\kappa:=\max\{1,\|1\|_{\mathcal A}\}$. For $n\geq1$, the power
$u\mapsto u^n$ is a continuous polynomial constructed from
multilinear multiplication. By
Proposition~\ref{prop:diferenciabilidad-operadores-multilineales},
for $0\leq j\leq n$ its derivative of order $j$ is
\[
D^j(u\mapsto u^n)(u)[h_1,\ldots,h_j]
=\frac{n!}{(n-j)!}u^{n-j}h_1\cdots h_j.
\]
For $j>n$ it is zero. If $\|u\|_{\mathcal A}\leq r$, the norm of
this derivative is at most
$\displaystyle\kappa\displaystyle\frac{n!}{(n-j)!}r^{n-j}$; the factor $\kappa$
also covers $n=j$, when the unit occurs.

Fix $0<r<\rho<R$ and $j\in\mathbb N_0$. The sequence
$\displaystyle\frac{n!}{(n-j)!}\left(\displaystyle\frac{r}{\rho}\right)^{n-j}$,
for $n\geq j$, is bounded. Indeed, the ratio of consecutive terms
tends to $\displaystyle\frac{r}{\rho}<1$; after some index it is
less than $\displaystyle\frac{1}{2}\left(1+\displaystyle\frac{r}{\rho}\right)<1$,
and there are only finitely many earlier terms. Since
$\displaystyle\sum_n|a_n|\rho^n<\infty$, it follows that
\[
\sum_{n=j}^\infty|a_n|\frac{n!}{(n-j)!}r^{n-j}<\infty.
\]
In particular, the error in truncating the series of derivatives after
$N$ is bounded, uniformly on the ball of radius $r$, by
\[
\kappa\sum_{n>N,\ n\geq j}
|a_n|\frac{n!}{(n-j)!}r^{n-j},
\]
which tends to zero as the tail of a convergent series of nonnegative
numbers.

Let $F_j$ be the uniform limit of the derivatives of order $j$ of the
partial polynomials. It is continuous as a uniform limit of continuous
maps. For $u,u+h$ in a smaller ball, the integral mean value formula
applied to the derivatives of order $j$ of each polynomial passes to
the limit by the preceding uniform bound and gives
\[
F_j(u+h)-F_j(u)
=\int_0^1 F_{j+1}(u+th)[h,\,\cdot\,]\,dt.
\]
The integral takes values in the space of $j$--linear operators;
for $j=0$ it takes values in $\mathcal A$. It is the integral of a
continuous curve with values in a Banach space. For $h\neq0$,
\[
\frac{\|F_j(u+h)-F_j(u)-F_{j+1}(u)[h,\,\cdot\,]\|}{\|h\|_{\mathcal A}}
\leq\sup_{0\leq t\leq1}\|F_{j+1}(u+th)-F_{j+1}(u)\|.
\]
Continuity at $u$ makes the right-hand side tend to zero: the entire
segment lies at distance at most $\|h\|_{\mathcal A}$ from $u$.
Therefore $DF_j=F_{j+1}$. Induction on $j$ proves smoothness and
term-by-term differentiation.
\end{proof}

\begin{example}[Elementary functions as maps between function spaces]
\label{ej:funciones-elementales-algebras}
Let $\mathcal A$ be one of the spaces in
Lemmas~\ref{lem:algebra-ck-compacta} and~\ref{lem:algebra-h1-intervalo},
with the submultiplicative norm fixed there. All its elements are
regarded as continuous functions on the corresponding compact set.
Then
\[
\mathcal A_+:=\{f\in\mathcal A\mid f(x)>0\text{ for every }x\}
\]
is open. The maps
\[
\operatorname{Exp}\colon\mathcal A\to\mathcal A_+,\quad f\mapsto e^f,
\qquad
\operatorname{Log}\colon\mathcal A_+\to\mathcal A,\quad f\mapsto\log f,
\]
are smooth and inverse to each other. On $\mathcal A_+$, the maps
$f\mapsto f^{-1}$ and $f\mapsto f^{-\displaystyle\frac{1}{2}}$ are also smooth.
In particular,
\[
\begin{aligned}
D\operatorname{Exp}(f)h&=e^fh,&
D\operatorname{Log}(f)h&=f^{-1}h,\\
D(f\mapsto f^{-1})(f)h&=-f^{-2}h,&
D(f\mapsto f^{-\frac{1}{2}})(f)h&=-\frac{1}{2}f^{-\frac{3}{2}}h.
\end{aligned}
\]
\end{example}

\begin{proof}
If $f\in\mathcal A_+$, its minimum $m_f$ on the compact set is
positive. The inequality $\|h\|_\infty\leq\|h\|_{\mathcal A}$ shows
that $f+h>0$ whenever
$\|h\|_{\mathcal A}<\displaystyle\frac{m_f}{2}$. This proves
openness. The pointwise functions $f^{-1}$, $\log f$, and
$f^{-\displaystyle\frac{1}{2}}$ belong to $\mathcal A$. For $C^k$, apply the
chain rule to smooth scalar functions on a neighborhood of the compact
interval $[m_f,\|f\|_\infty]$. For $H^1$, take the smooth approximations
$f_\nu$ constructed in Lemma~\ref{lem:algebra-h1-intervalo}.
They converge uniformly to $f$, and their derivatives converge to $f'$
in $L^2$. For sufficiently large $\nu$, their values lie in the
compact interval $[\displaystyle\frac{m_f}{2},\|f\|_\infty+1]$.
For any of the three scalar functions $\psi$ under consideration,
$\psi(f_\nu)$ converges uniformly to $\psi(f)$ and
\[
\|\psi'(f_\nu)f_\nu'-\psi'(f)f'\|_{L^2}
\leq\|\psi'(f_\nu)\|_\infty\|f_\nu'-f'\|_{L^2}
+\|\psi'(f_\nu)-\psi'(f)\|_\infty\|f'\|_{L^2}
\longrightarrow0.
\]
The first uniform norm is bounded by the maximum of $|\psi'|$ on
that interval, and the second tends to zero by uniform continuity
of $\psi'$. Passing to the limit against test functions identifies
the weak derivative of $\psi(f)$ as $\psi'(f)f'$.
Thus the three weak derivatives are, respectively,
\[
-f^{-2}f',\qquad f^{-1}f',\qquad
-\frac{1}{2}f^{-\frac{3}{2}}f'.
\]
They belong to $L^2$, since the factors multiplying $f'$ are bounded.
Thus all three compositions belong to $H^1$.

Lemma~\ref{lem:series-potencias-algebra-banach} applies to the series
\[
\begin{aligned}
\operatorname{Exp}(u)&=\sum_{n=0}^\infty\frac{u^n}{n!},\\
L(w)&=\sum_{n=1}^\infty\frac{(-1)^{n+1}}{n}w^n,
&\|w\|_{\mathcal A}&<1,\\
I(w)&=\sum_{n=0}^\infty(-1)^nw^n,
&\|w\|_{\mathcal A}&<1,\\
S(w)&=\sum_{n=0}^\infty
(-1)^n\frac{\binom{2n}{n}}{4^n}w^n,
&\|w\|_{\mathcal A}&<1.
\end{aligned}
\]
The first has infinite radius, and the others have radius one.
Evaluation is continuous, so the scalar identities for these series
give $L(w)=\log(1+w)$, $I(w)=(1+w)^{-1}$, and
$S(w)=(1+w)^{-\displaystyle\frac{1}{2}}$ pointwise. These scalar identities can
be checked by differentiating the series for $|t|<1$: that for $L$
has derivative $\displaystyle\frac{1}{1+t}$ and value zero at the
origin; that for $I$, multiplied by $1+t$, equals one; and that for
$S$ satisfies $2(1+t)S'(t)+S(t)=0$, $S(0)=1$. The last equation
implies that $\sqrt{1+t}S(t)$ is constant and hence equals one.

For a fixed $f\in\mathcal A_+$, if
$\|h\|_{\mathcal A}<\displaystyle\frac{1}{2(1+\|f^{-1}\|_{\mathcal A})}$,
then $\|hf^{-1}\|_{\mathcal A}<\displaystyle\frac{1}{2}$ and
\[
\begin{aligned}
\log(f+h)&=\log f+L(hf^{-1}),\\
(f+h)^{-1}&=f^{-1}I(hf^{-1}),\\
(f+h)^{-\frac{1}{2}}&=f^{-\frac{1}{2}}S(hf^{-1}).
\end{aligned}
\]
These are compositions of smooth maps between Banach spaces.
The local formulas prove smoothness throughout $\mathcal A_+$ and
give the stated derivatives. The exponential series yields, for
$j\geq1$,
\[
D^j\operatorname{Exp}(f)[h_1,\ldots,h_j]
=e^fh_1\cdots h_j.
\]
Finally, $\log(e^f)=f$ and $e^{\log f}=f$ pointwise; uniqueness of
the continuous representatives makes these identities in $\mathcal A$.
\end{proof}
\section{Lagrange multipliers}
\label{sec:multiplicadores-lagrange-calculo-banach}

Theorem~\ref{teo:funcion-implicita-banach} describes a regular level set as a graph over the kernel of the derivative of the constraint. In these coordinates, a constrained extremum becomes an ordinary extremum on the tangent space. This observation leads to the principle of Lagrange multipliers and explains why the derivative of the functional must be a combination of the derivatives of the constraints.

\begin{theorem}[Lagrange multipliers for finitely many constraints]
\label{teo:multiplicadores-lagrange-banach-general}
\index{Lagrange multipliers}
Let $X$ be a Banach space, let $U\subseteq X$ be open, and let $J\in C^1(U,\mathbb R)$ and $F=(F_1,\ldots,F_m)\in C^1(U,\mathbb R^m)$. Let $x_0\in U$ be a local extremum of $J\restriction_{F^{-1}(F(x_0))}$ and suppose $DF(x_0)\colon X\longrightarrow\mathbb R^m$ is surjective. Then there exists a unique functional $\lambda\in(\mathbb R^m)'$ such that
\begin{equation}
\label{eq:multiplicadores-lagrange-finito}
DJ(x_0)=\lambda\circ DF(x_0).
\end{equation}
Equivalently, there exist unique numbers $\lambda_1,\ldots,\lambda_m\in\mathbb R$ such that
\[
DJ(x_0)=\displaystyle\sum_{j=1}^{m}\lambda_jDF_j(x_0).
\]
\end{theorem}

\begin{proof}
Let $K:=\ker DF(x_0)$. Since $DF(x_0)$ is surjective and $\mathbb R^m$ is finite-dimensional, for each vector $e_j$ in the standard basis of $\mathbb R^m$ there exists $v_j\in X$ such that $DF(x_0)v_j=e_j$. Define $E:=\operatorname{span}\{v_1,\ldots,v_m\}$. Then $E$ is finite-dimensional and therefore closed. Moreover, $DF(x_0)\restriction_E\colon E\longrightarrow\mathbb R^m$ is an isomorphism. For each $v\in X$, if
\[
w:=\bigl(DF(x_0)\restriction_E\bigr)^{-1}(DF(x_0)v),
\]
then $v-w\in K$. Thus $X=K\oplus E$.

Consider the open set
\[
\mathcal O:=\{(u,v)\in K\times E\mid x_0+u+v\in U\}
\]
and the map $G\colon \mathcal O\longrightarrow\mathbb R^m$ defined by
$G(u,v):=F(x_0+u+v)-F(x_0)$. Moreover, $G(0,0)=0$ and
\[
D_2G(0,0)=DF(x_0)\restriction_E\colon E\longrightarrow\mathbb R^m.
\]
Theorem~\ref{teo:funcion-implicita-banach} gives a neighborhood $W\subseteq K$ of $0$ and a function $\psi\colon W\longrightarrow E$ of class $C^1$ such that $\psi(0)=0$ and
\[
F(x_0+u+\psi(u))=F(x_0)
\]
for every $u\in W$. Differentiating this identity at $u=0$ gives
\[
DF(x_0)\bigl(u+D\psi(0)u\bigr)=0
\]
for every $u\in K$. Since $DF(x_0)u=0$ and $D\psi(0)u\in E$, injectivity of $DF(x_0)\restriction_E$ implies $D\psi(0)u=0$. Thus $D\psi(0)=0$.

Fix $u\in K$. For every $t$ in a neighborhood of $0$, the curve $\gamma(t):=x_0+tu+\psi(tu)$ lies in $F^{-1}(F(x_0))$. Since $x_0$ is a local extremum of the restriction of $J$, the real function $J\circ\gamma$ has a local extremum at $0$. Fermat's rule and the chain rule give
\[
0=(J\circ\gamma)'(0)=DJ(x_0)\bigl(u+D\psi(0)u\bigr)=DJ(x_0)u.
\]
Thus $K\subseteq\ker DJ(x_0)$.

Let $R\colon \mathbb R^m\longrightarrow E$ be the inverse of $DF(x_0)\restriction_E$ and define $\lambda:=DJ(x_0)\circ R\in(\mathbb R^m)'$. For each $v\in X$, the vector $v-R(DF(x_0)v)$ belongs to $K$; hence
\[
DJ(x_0)v
=
DJ(x_0)R(DF(x_0)v)
=
\lambda(DF(x_0)v).
\]
This proves \eqref{eq:multiplicadores-lagrange-finito}. If $\mu\in(\mathbb R^m)'$ also satisfies $DJ(x_0)=\mu\circ DF(x_0)$, surjectivity of $DF(x_0)$ implies $\lambda=\mu$. The formulation in terms of the numbers $\lambda_j$ follows by identifying $(\mathbb R^m)'$ with $\mathbb R^m$ using the standard basis.
\end{proof}

\begin{corollary}[Lagrange multipliers for one constraint]
\label{teo:multiplicadores-lagrange-banach}
\index{method!of Lagrange multipliers}
Let $X$ be a Banach space, let $U\subseteq X$ be open, and let $J,F\in C^1(U,\mathbb R)$. If $x_0$ is a local extremum of $J$ restricted to $\{x\in U\mid F(x)=F(x_0)\}$ and $DF(x_0)\neq0$, there exists a unique $\lambda\in\mathbb R$ such that $DJ(x_0)=\lambda DF(x_0)$.
\end{corollary}

\begin{proof}
The operator $DF(x_0)\colon X\longrightarrow\mathbb R$ is surjective because it is linear and nonzero. Apply the preceding theorem with $m=1$ and identify each linear functional on $\mathbb R$ with multiplication by a unique real number.
\end{proof}

\subsection{Exercises}
\label{subsec:ejercicios-teoremas-fundamentales-calculo}

\begin{exercise}
\label{ejer:inversa-sin-continuidad-derivada}
Define $f\colon \mathbb R\longrightarrow\mathbb R$ by $f(0)=0$ and $\displaystyle f(x)=x+2x^2\sin\left(\displaystyle\frac{1}{x}\right)$ for $x\neq0$. Prove that $f'(0)=1$ but $f$ is not injective on any neighborhood of $0$. Determine which assumption of Theorem~\ref{teo:funcion-inversa-local-banach} fails.
\end{exercise}

\begin{exercise}
\label{ejer:coordenadas-polares-laplaciano}
Let $\Phi(r,\theta):=(r\cos\theta,r\sin\theta)$ and let $u\in C^2(\mathbb R^2\setminus\{0\})$. If $v:=u\circ\Phi$, use the second-order chain rule to prove that
\[
(\Delta u)\circ\Phi
=
\frac{\partial^2v}{\partial r^2}
+
\frac{1}{r}\frac{\partial v}{\partial r}
+
\frac{1}{r^2}\frac{\partial^2v}{\partial\theta^2}.
\]
Interpret this identity as a first example of the behavior of a differential operator under a change of coordinates.
\end{exercise}

\begin{exercise}
\label{ejer:dependencia-raiz-simple}
Let $P(z,a_0,\ldots,a_{n-1}):=z^n+a_{n-1}z^{n-1}+\cdots+a_0$. Suppose $z_0$ is a simple root of the polynomial corresponding to $(a_0^0,\ldots,a_{n-1}^0)$. Apply Theorem~\ref{teo:funcion-implicita-banach} to prove that this root depends smoothly on the coefficients in a neighborhood of the given point. Compute the derivative with respect to each coefficient.
\end{exercise}

\begin{exercise}
\label{ejer:esfera-hilbert-subvariedad-lagrange}
Let $H$ be a real Hilbert space and let $S:=\{u\in H\mid\|u\|_H=1\}$. Use Theorem~\ref{teo:funcion-implicita-banach} to prove that $S$ is a Banach submanifold of codimension one and that $T_uS=u^\perp$. If $J\in C^1(H,\mathbb R)$ has a local extremum on $S$, prove that there exists $\lambda\in\mathbb R$ such that $DJ(u)(v)=\lambda\langle u,v\rangle_H$ for every $v\in H$.
\end{exercise}

\begin{exercise}
\label{ejer:lagrange-rayleigh-haz}
Let $(M,\mathbf{g})$ be a compact Riemannian manifold with smooth boundary, and let $\mathbf{E}\to M$ be a smooth real vector bundle of finite rank, equipped with a bundle metric $\mathbf{h}_{\mathbf{E}}$ and a compatible connection $\nabla^{\mathbf{E}}$. On $W_0^{1,2}(M,\mathbf{E})$, consider the functionals
$\displaystyle J(\mathbf{u}):=\displaystyle\frac{1}{2}\int_M|\nabla_w^{\mathbf{E}}\mathbf{u}|_{\mathbf{g},\mathbf{h}_{\mathbf{E}}}^2\,d\lambda_{\mathbf{g}}$ and $\displaystyle G(\mathbf{u}):=\displaystyle\frac{1}{2}\int_M|\mathbf{u}|_{\mathbf{h}_{\mathbf{E}}}^2\,d\lambda_{\mathbf{g}}$. Suppose $\mathbf{u}\neq0$ minimizes $J$ subject to the constraint $G(\mathbf{u})=1$. Apply the Lagrange multiplier theorem to obtain the weak identity for the eigenvalue problem associated with the Bochner Laplacian.
\end{exercise}

\begin{exercise}
\label{ejer:implicita-ecuacion-secciones}
Let $X$ and $Y$ be Banach spaces of sections of vector bundles, and let $\mathcal F\colon X\times\mathbb R\longrightarrow Y$ be of class $C^1$. Suppose $\mathcal F(\mathbf{u}_0,\lambda_0)=0$ and $D_1\mathcal F(\mathbf{u}_0,\lambda_0)\colon X\longrightarrow Y$ is an isomorphism. Use Theorem~\ref{teo:funcion-implicita-banach} to formulate and prove existence of a local branch $\lambda\longmapsto \mathbf{u}(\lambda)$ of solutions. Compute $\mathbf{u}'(\lambda_0)$.
\end{exercise}

\section{Fréchet spaces and calculus by directional derivatives}
\label{sec:calculo-frechet-autocontenido}

The Fréchet derivative studied so far is defined by means of a norm.
The name has a different use in the expression \emph{Fréchet space}: it
refers to locally convex spaces that need not be normable. To work, for
example, with all derivatives of a smooth function, we must specify both
the topology of the space and the calculus used on it. This does not change
the definitions of class $C^k$ already fixed for Banach spaces.

Moore's exposition \cite[Sections 1.2 and 1.3]{Moore2017GlobalAnalysis}
motivates this step through chains of Banach spaces. We use the calculus
of jointly continuous directional derivatives, usually called Bastiani
calculus and also used by Hamilton \cite[Part I]{Hamilton1982NashMoser}.

\subsection{Conventions and models}
\label{subsec:seminormas-frechet}

We use the notion of a Fréchet space from
Definition~\ref{def: espacio de Frechet}, its seminorm characterization
(Theorem~\ref{teo: caracterizacion Frechet por seminormas}), and the
convergence and continuity criteria in Section~\ref{ap:espacios-de-Frechet}
of Appendix~\ref{ap:analisis-matematico} and Propositions~\ref{prop: convergencia y Cauchy por seminormas} and
\ref{prop: metrizacion por seminormas}. In each space, fix a countable
separating family of seminorms $(p_j)_{j\in\mathbb N_0}$.
We may take it increasing by replacing $p_j$ with
$\displaystyle\max_{0\leq i\leq j}p_i$; each family controls the other
through finite maxima. The index controlling a seminorm and its associated
constant may depend on that seminorm, as specified in
Corollary~\ref{prop:metrica-continuidad-frechet}.

The test function spaces $\mathcal D_K(\Omega)$ and $\mathcal S(\mathbb R^n)$
were studied in Chapter~\ref{cap:distribuciones-fourier}; for the latter, see
Theorem~\ref{teo:schwartz-frechet}. The geometric applications also require
the following models.

\begin{example}[Products and smooth sections]

\label{ej:modelos-frechet}

Every Banach space is Fréchet, by taking $p_j=\|\cdot\|$.
The product $\mathbb R^{\mathbb N_0}$ is Fréchet with
$p_j(x)=\max_{0\leq i\leq j}|x_i|$: a Cauchy sequence converges in each
coordinate, and these limits form an element of the product.

Let $\mathbf E\to M$ be a smooth finite-rank vector bundle over a
second-countable finite-dimensional smooth manifold. On $\Gamma(\mathbf E)$,
take the seminorms controlling derivatives of the components up to a given
order on coordinate compact sets. A countable family suffices: choose a
compact exhaustion and, at each stage, finitely many relatively compact
charts and trivializations covering that compact set. Taking cumulative
maxima gives increasing seminorms. The resulting topology is uniform
convergence on compact sets of all derivatives.

Proposition~\ref{prop: GammaK Frechet} proves completeness for smooth
sections supported in a fixed compact set. To apply it to sections
without a support restriction, choose a smooth locally finite partition
of unity $(\chi_i)_{i\in\mathbb N}$ with compact
$K_i:=\operatorname{supp}\chi_i$. If $(\mathbf u_n)$ is Cauchy in the
topology of convergence on compact sets, the Leibniz rule shows that
$(\chi_i\mathbf u_n)_n$ is Cauchy in $\Gamma_{K_i}(\mathbf E)$ for each
$i$. Let $\mathbf v_i$ be its limit. The sum
\[
 \mathbf u=\sum_{i\in\mathbb N}\mathbf v_i
\]
is locally finite and defines a smooth section. Only finitely many
indices $i$ occur on a compact set $K$; the identity
$\mathbf u_n=\displaystyle \sum_{i\in\mathbb N}\chi_i\mathbf u_n$ then gives convergence
to $\mathbf u$ with every derivative on $K$. This proves completeness.
For a compact base, take $K=M$ directly in the cited proposition.
Equivalence between smoothness up to the boundary and local smooth
extension is proved in Lemma~\ref{lem:prolongacion-suave-frontera-frechet}.

If $M$ is compact, we may equivalently use

\[
 p_j(\mathbf u)=\displaystyle\sum_{i=0}^{j}
      \displaystyle\sup_{x\in M}|\nabla^i\mathbf u(x)|,
 \qquad j\in\mathbb N_0.
\]

Theorem~\ref{norma ck independiente de metrica y conexion}, applied to
the compactly supported cutoffs, proves independence of the connection,
metric and trivialization for this topology. The comparisons hold on
each compact set; they do not assert equivalence of global norms on a noncompact base.

\end{example}

\begin{lemma}[Local smooth extension at the boundary]

\label{lem:prolongacion-suave-frontera-frechet}

Let $u(x',t)$ be defined on a half-ball $t\geq0$, with partial
derivatives of every order in the interior extending continuously to the
boundary. Near each point of $t=0$, $u$ is the restriction of a smooth
function defined on an open subset of $\mathbb R^n$. The assertion applies
componentwise to finite-rank sections.

\end{lemma}

\begin{proof}

Work in a smaller cylinder. Integral identities in tangential directions,
passed to the limit from $t>0$, show that $a_j(x')=\partial_t^ju(x',0)$
is smooth and its tangential derivatives are traces of the mixed
derivatives of $u$. Multiplying all $a_j$ by a tangential cutoff equal
to one near the chosen point gives supports in one compact set.
Let $\chi\in C_c^\infty(\mathbb R)$ equal one near zero. We construct
a function with these jets:

\[
 S(x',t)=\displaystyle\sum_{j=0}^{\infty}
             a_j(x')\frac{t^j}{j!}\chi(t/\varepsilon_j).
\]

For $j\geq1$, choose $0<\varepsilon_j<1/j$ small enough that the
$C^m$ norm of the $j$th summand is at most $2^{-j}$ for every
$m\leq\lfloor j/2\rfloor$. This is possible: a total derivative of
order at most $m$ contains factors bounded by constants times
$\varepsilon_j^{j-b}$, where $b\leq m$ is the number of normal
derivatives. Since $j-m>0$, all these terms tend to zero with
$\varepsilon_j$. There are only finitely many requirements for each
$j$. Choose any $\varepsilon_0>0$ for the zeroth term.

For each fixed $m$, the tail with $j\geq2m+1$ converges in $C^m$ by
the Weierstrass test. Thus $S$ is smooth and the series can be
differentiated termwise. Each cutoff is constant near zero, so
$\partial_t^kS(x',0)=a_k(x')$ for every $k$. Define the extension to
be $u$ for $t\geq0$ and $S$ for $t<0$, in the cylinder where the
tangential cutoff equals one. All partial derivatives agree on $t=0$.
The fundamental theorem of calculus in the normal direction first shows
that the glued function is $C^1$; applied successively to its partial
derivatives, it proves $C^m$ regularity for every $m$. This gives the
required smooth extension, given by a single series realizing all
prescribed jets.

\end{proof}

\begin{proposition}[A topology not induced by a norm]

\label{prop:no-normabilidad-Cinfty}

The preceding topology on $C^\infty([0,1],\mathbb R)$ is not induced
by any norm.

\end{proposition}

\begin{proof}

Suppose that a norm $\|\cdot\|_*$ induces it. Continuity of the
identity gives $m$ and $\delta>0$ such that $p_m(u)<\delta$ implies
$\|u\|_*<1$. There is also $C>0$ with $p_{m+1}(u)\leq C\|u\|_*$.
For integers $n\geq2$, set

\[
 u_n(t)=\frac{\delta}{2(m+1)}(2\pi n)^{-m}\sin(2\pi nt).
\]

Then $p_m(u_n)\leq\delta/2$, whereas
$p_{m+1}(u_n)\geq\delta(2\pi n)/(2(m+1))\to\infty$.
This contradicts $p_{m+1}(u_n)\leq C$. Thus the issue is not that we
have chosen an incomplete norm: no norm produces this topology.

\end{proof}

\subsection{Integration and differential calculus}
\label{subsec:integracion-calculo-bastiani}

\begin{lemma}[Integration of Fréchet-valued curves]

\label{lem:integracion-curvas-frechet}

If $a\colon[\alpha,\beta]\to F$ is continuous, its Riemann sums
converge in $F$ to an element denoted by $\displaystyle \int_\alpha^\beta a(t)\,dt$.
For every continuous seminorm $q$ and every continuous linear operator
$A\colon F\to G$, one has

\[
 q\left(\int_\alpha^\beta a(t)\,dt\right)
 \leq\int_\alpha^\beta q(a(t))\,dt,
 \qquad
 A\left(\int_\alpha^\beta a(t)\,dt\right)
 =\int_\alpha^\beta Aa(t)\,dt.
\]

The function $b(t)=\displaystyle \int_\alpha^t a(s)\,ds$ is continuously
differentiable and $b'(t)=a(t)$. If $a$ is continuously differentiable,
then $a(\beta)-a(\alpha)=\displaystyle \int_\alpha^\beta a'(t)\,dt$.

\end{lemma}

\begin{proof}

Fix a seminorm $q$. Continuity on the compact interval is uniform
for the pseudodistance $q(a(t)-a(s))$. Comparing a Riemann sum with the
sum for a refinement bounds their difference by
$(\beta-\alpha)\omega_q(\delta)$, where $\omega_q(\delta)\to0$ as the
mesh $\delta$ tends to zero. Two partitions are compared through a common
refinement. A sequence of sums with meshes tending to zero is therefore
Cauchy for each seminorm and converges by completeness. The same estimate
shows that the limit is independent of partitions and sample points.
The two assertions in the display follow by passing to the limit in the sums.

Moreover,

\[
 q\left(\frac{b(t+h)-b(t)}{h}-a(t)\right)
 \leq\displaystyle\sup_{|s-t|\leq|h|}q(a(s)-a(t))\longrightarrow0.
\]

This proves differentiation of the integral. For the final identity,
the difference of its two sides is annihilated by every continuous linear
functional, by the real fundamental theorem of calculus. Such functionals
separate points: for $v\neq0$, choose a seminorm $q$ with $q(v)>0$, define
$\ell(tv)=tq(v)$, and apply Hahn--Banach with bound $q$.
Separation proves the identity in $F$.

\end{proof}

\begin{definition}[Bastiani differentiability]

\label{def:calculo-bastiani}

Let $E,F$ be Fréchet spaces and let $U\subseteq E$ be open. A map
$f\colon U\to F$ is of class $C^1_{\mathrm B}$ if it is continuous,

\[
 Df(x)v=\lim_{t\to0}\frac{f(x+tv)-f(x)}{t},
 \qquad x\in U,\quad v\in E,
\]

exists for every $(x,v)$, with the limit taken in $F$, and
$(x,v)\mapsto Df(x)v$ is jointly continuous. For $r\geq2$, require the
iterated derivatives

\[
 D^j f(x)(v_1,\ldots,v_j)
 =D_{v_j}\bigl(D^{j-1}f(\,\cdot\,)(v_1,\ldots,v_{j-1})\bigr)(x),
 \quad 1\leq j\leq r,
\]

to exist and be jointly continuous on $U\times E^j$.
Write $f\in C^\infty_{\mathrm B}(U,F)$ if these conditions hold for every
$r$. In the Fréchet category, \emph{smooth} means $C^\infty_{\mathrm B}$.
No norm is imposed on the operator space $\mathcal L(E,F)$.

\end{definition}

\begin{proposition}[Linearity and the integral formula]

\label{prop:linealidad-bastiani}

If $f\in C^1_{\mathrm B}(U,F)$, then $Df(x)\colon E\to F$ is linear
and continuous. Whenever $x+[0,1]h\subseteq U$,

\begin{equation}\label{eq:valor-medio-frechet}
 f(x+h)-f(x)=\int_0^1Df(x+th)h\,dt.
\end{equation}

\end{proposition}

\begin{proof}

The function $t\mapsto f(x+th)$ has continuous derivative
$Df(x+th)h$ by definition. The preceding lemma proves the integral formula
without yet using linearity in $h$. Homogeneity of $Df(x)$ follows by a
change of variable in the limit. For additivity, split

\begin{align*}
 \frac{f(x+t(v+w))-f(x)}{t}
 &=\frac{f(x+tv+tw)-f(x+tv)}{t}
   +\frac{f(x+tv)-f(x)}{t}\\
 &=\int_0^1Df(x+tv+stw)w\,ds
   +\frac{f(x+tv)-f(x)}{t}.
\end{align*}

For small $t$, the segments used lie in a convex neighborhood of $x$
contained in $U$. In each seminorm, joint continuity makes the integrand
converge to $Df(x)w$ uniformly for $s\in[0,1]$. Taking limits gives
$Df(x)(v+w)=Df(x)v+Df(x)w$. Continuity of this operator follows by
restricting the continuous map $(x,v)\mapsto Df(x)v$ to a fixed $x$.

\end{proof}

\begin{theorem}[Rules of Fréchet calculus]

\label{teo:reglas-calculo-frechet}

Compositions of $C^r_{\mathrm B}$ maps are $C^r_{\mathrm B}$, and
their differentials satisfy the chain rule. The derivatives $D^j f(x)$
are continuous symmetric multilinear maps. If
$B\colon F_1\times F_2\to G$ is continuous bilinear and $f_i\colon U\to F_i$
are $C^r_{\mathrm B}$, then $B(f_1,f_2)$ is $C^r_{\mathrm B}$ and

\[
 D(B(f_1,f_2))(x)v
 =B(Df_1(x)v,f_2(x))+B(f_1(x),Df_2(x)v).
\]

Furthermore, if $f$ is smooth, all its derivative evaluation maps
$(x,v_1,\ldots,v_j)\mapsto D^jf(x)(v_1,\ldots,v_j)$ are smooth.

\end{theorem}

\begin{proof}

Let $f\colon U\to V\subseteq F$ and $g\colon V\to G$ be
$C^1_{\mathrm B}$. Put $a_t=f(x+tv)-f(x)$. For small $t$,
$f(x)+[0,1]a_t\subseteq V$, and the integral formula gives

\[
 \frac{g(f(x+tv))-g(f(x))}{t}
 =\int_0^1Dg(f(x)+sa_t)\frac{a_t}{t}\,ds
 \longrightarrow Dg(f(x))Df(x)v.
\]

The limit is taken in each seminorm of $G$ and is uniform in $s$:
$a_t\to0$ and $a_t/t\to Df(x)v$, so the arguments of the integrand
converge uniformly to the same point. Joint continuity controls them in
a neighborhood of that point. The resulting differential is jointly
continuous in $(x,v)$.

The bilinear rule follows from the exact identity
$B(a+\Delta a,b+\Delta b)-B(a,b)
=B(\Delta a,b)+B(a,\Delta b)+B(\Delta a,\Delta b)$.
The last term divided by $t$ tends to zero. To justify this by seminorms,
bilinear continuity gives, for every target seminorm $q$, continuous
seminorms $p_1,p_2$ and $C>0$ such that
$q(B(a,b))\leq Cp_1(a)p_2(b)$. This estimate follows by scaling two zero neighborhoods, as in Theorem~\ref{teo: continuidad mediante seminormas}.

Linearity in the last direction of an iterated derivative follows by
applying the preceding proposition to the previous derivative with the
other directions fixed. For symmetry, restrict to
$x+t_1v_1+\cdots+t_jv_j$ on a small open set in $\mathbb R^j$.
After composing with any continuous linear functional on the target,
continuity of the partial derivatives and the real Schwarz theorem allow
them to be permuted. Continuous linear functionals separate points by the
proof of the integration lemma. This gives symmetry in $F$ and hence
linearity in every direction. Multilinear continuity is already part of
the joint continuity required in the definition.

To complete the higher-order assertions, no topology on
$\mathcal L(E,F)$ is needed. The differential of the evaluation map
$J_jf(x,v_1,\ldots,v_j)=D^jf(x)(v_1,\ldots,v_j)$ is

\begin{align*}
 D(J_jf)(x,v_1,\ldots,v_j)(h,w_1,\ldots,w_j)
 &=D^{j+1}f(x)(v_1,\ldots,v_j,h)\\
 &\quad+\displaystyle\sum_{i=1}^j
 D^jf(x)(v_1,\ldots,v_{i-1},w_i,v_{i+1},\ldots,v_j).
\end{align*}

Indeed, first vary $x$ and then each direction. Multilinearity shows
that terms with two direction increments have order at least two in the
parameter, and joint continuity controls their coefficients. This formula
proves smoothness of $J_jf$ inductively. Differentiating the chain rule
produces finite sums of compositions of such evaluation maps. Induction
proves the composition rule at every stated order.

\end{proof}

\begin{proposition}[Taylor formula with integral remainder]

\label{prop:taylor-frechet}

If $f\in C^r_{\mathrm B}(U,F)$, $r\geq1$, and $x+[0,1]h\subseteq U$,
then

\begin{align*}
 f(x+h)&=\displaystyle\sum_{j=0}^{r-1}\frac1{j!}D^jf(x)(h,\ldots,h)\\
 &\quad+\frac1{(r-1)!}\int_0^1(1-t)^{r-1}
                  D^rf(x+th)(h,\ldots,h)\,dt.
\end{align*}

Here $D^0f(x)=f(x)$. If the order $r$ term at $x$ is also subtracted,
the remainder $R_r(x,h)$ satisfies, for every continuous seminorm $q$,

\[
 q(R_r(x,h))\leq\frac1{r!}\displaystyle\sup_{0\leq t\leq1}
 q\bigl((D^rf(x+th)-D^rf(x))(h,\ldots,h)\bigr).
\]

\end{proposition}

\begin{proof}

Apply the fundamental theorem of calculus successively to the curve
$a(t)=f(x+th)$. Its derivatives are
$a^{(j)}(t)=D^jf(x+th)(h,\ldots,h)$ by the chain rule.
Integrating $r$ times over the ordered parameter simplex produces the
factor $(1-t)^{r-1}/(r-1)!$. The identity may first be checked after
composition with each continuous linear functional; separation of points
then gives the identity in $F$. The estimate follows from the integration
lemma and $\displaystyle \int_0^1(1-t)^{r-1}\,dt=1/r$.
In particular, given $x$, a direction $v$, a continuous seminorm $q$,
and $\varepsilon>0$, there is $\delta>0$ such that $x+tv\in U$ and
\[
 q(R_r(x,tv))\leq\varepsilon |t|^r,\qquad |t|<\delta.
\]
Indeed, continuity of $z\mapsto D^rf(z)(v,\ldots,v)$ bounds the
integrand difference by $r!\varepsilon$ on the entire segment
$x+[0,1](tv)$ after reducing $\delta$. This radius may depend on
$x,v,q$, and $\varepsilon$.

\end{proof}

\subsection{Relation to Banach calculus and limits of the theory}
\label{subsec:comparacion-bastiani-banach}

\begin{proposition}[Agreement in the smooth class, not in class $C^1$]

\label{prop:comparacion-calculos-banach-frechet}

On open subsets of Banach spaces, every $C^r$ map in the sense of
Definition~\ref{def:funcion-Ck-banach} is $C^r_{\mathrm B}$.
Conversely, $C^{r+1}_{\mathrm B}$ implies $C^r$ for $r\geq1$.
Consequently, the two notions of smoothness agree.
The condition $C^1_{\mathrm B}$ alone does not imply norm Fréchet
differentiability.

\end{proposition}

\begin{proof}

The first implication follows by evaluating continuous multilinear
operators on their arguments. For the converse, joint continuity of
$D^{j+1}f$ at $(x_0,0,\ldots,0)$ and multilinearity give a convex
neighborhood $V$ of $x_0$ and a constant $C_j$ such that

\[
 \|D^{j+1}f(x)(v_1,\ldots,v_{j+1})\|
 \leq C_j\prod_{i=1}^{j+1}\|v_i\|,\qquad x\in V.
\]

Obtain this bound by scaling each direction into a ball where joint
continuity gives a bound by $1$. The integral formula, applied with the
directions fixed, implies

\[
 \|D^jf(y)-D^jf(x)\|_{\mathcal L^j}
 \leq C_j\|y-x\|,\qquad x,y\in V,
 \quad 1\leq j\leq r.
\]

Intersecting the neighborhoods corresponding to $j\in\{1,\ldots,r\}$
gives a common convex open set $V$ for all these bounds. Thus these
derivatives are continuous in operator norm. For $x,x+h\in V$ and $1\leq j\leq r$, the integral formula and the preceding bound give
\[
 \bigl\|D^{j-1}f(x+h)-D^{j-1}f(x)-D^jf(x)[h,\,\cdot\,]\bigr\|_{\mathcal L^{j-1}}
 \leq \frac{C_j}{2}\|h\|^2.
\]
For $j=1$, interpret $\mathcal L^0$ as the target space. Convexity of $V$
keeps the segment within the domain of the bound. Dividing by $\|h\|$
and letting $h\to0$ successively identifies the first $r$ Fréchet
derivatives. Hence $f$ is $C^r$ in the normed-space sense.

For the final assertion, consider $H=\ell^2(\mathbb N;\mathbb R)$,
$\varphi(t)=t^2/(1+t^2)$, and

\[
 f(x)_n=\frac1n\varphi(nx_n),\qquad
 (Df(x)v)_n=\varphi'(nx_n)v_n.
\]

Since $\varphi(0)=0$ and $L:=\displaystyle \sup_{t\in\mathbb R}|\varphi'(t)|<\infty$,
$f$ is well defined and $L$-Lipschitz. The difference quotient converges
to the displayed expression in $\ell^2(\mathbb N)$ by dominated convergence, with
$L^2|v_n|^2$ as a bound for the squares. If $x^{(a)}\to x$ and
$v^{(a)}\to v$ in $\ell^2(\mathbb N)$, the difference of derivative values is bounded
by $L\|v^{(a)}-v\|_2$ plus

\[
 \left(\displaystyle\sum_{n=1}^{\infty}
 |\varphi'(nx_n^{(a)})-\varphi'(nx_n)|^2|v_n|^2\right)^{1/2},
\]

which tends to zero by dominated convergence. Hence $f$ is
$C^1_{\mathrm B}$ and $Df(0)=0$. However, for $h_n=e_n/n$,

\[
 \frac{\|f(h_n)-f(0)-Df(0)h_n\|_2}{\|h_n\|_2}
 =\varphi(1)=\frac12.
\]

Thus no norm Fréchet derivative exists at the origin. The continuity
required in Proposition~\ref{prop:continuidad-gateaux-implica-frechet}
is operator-norm continuity; it cannot be replaced by joint continuity
of the derivative evaluation map.

\end{proof}

\begin{example}[Why inverse and flow theorems do not transfer]

\label{ej:fallos-inversa-flujo-frechet}

On $E=\mathbb R^{\mathbb N_0}$ with the product topology, the map
$F(x)_j=e^{x_j}$ is smooth. Each differential
$DF(x)v=(e^{x_j}v_j)_j$ is a linear topological isomorphism, with inverse
$w\mapsto(e^{-x_j}w_j)_j$. These assertions are checked on each finite
set of coordinates, which are all that a product seminorm controls.
However, $F(E)$ contains no nonempty open subset of $E$: a basic
neighborhood restricts only finitely many coordinates, and one of the
remaining coordinates can be made negative. Since $F(E)$ consists of
strictly positive sequences, $F$ is not a local diffeomorphism.
Even invertibility of the differential at every point is insufficient.

Existence and uniqueness can also fail for smooth vector fields.
The field $X(x)_j=x_j^2$, with initial data $x_j(0)=j+1$, would force
$x_j(t)=(j+1)/(1-(j+1)t)$ on the existence interval of each coordinate.
There is no common interval $[0,\varepsilon)$ with $\varepsilon>0$ for
all coordinates. For uniqueness, take the continuous linear operator
$S(x)_j=x_{j+1}$ and a nonzero smooth function $a\colon\mathbb R\to\mathbb R$
identically zero for $t\leq0$, for instance $a(t)=e^{-1/t^2}$ for $t>0$.
Then

\[
 u(t)=(a^{(j)}(t))_{j\in\mathbb N_0},\qquad
 u'(t)=Su(t),\qquad u(0)=0,
\]

is a smooth solution distinct from the zero solution. A curve with
values in the product is smooth exactly when every coordinate is smooth,
as follows from the definition and the finite-coordinate seminorms.
These examples explain why the results in the following chapters are
applied at a Banach or Hilbert level when inversion, differential equations,
or descent are involved.

\end{example}

\chapter{Banach and Fréchet manifolds}
\label{cap:variedades-banach}

Infinite-dimensional manifolds extend the notion of a smooth manifold beyond the finite-dimensional setting. Whereas a finite-dimensional manifold is described locally by open subsets of $\mathbb{R}^n$, in infinite dimensions various topological vector spaces can serve as local models. We first develop Banach theory, keeping complementability and regularity hypotheses explicit. Section~\ref{sec:variedades-frechet-autocontenidas} then constructs the smooth Fréchet category, specifying which local results survive and which do not. ILH towers and their mapping spaces are studied in Chapter~\ref{cap:geometria-uniforme-funcional}.

This extension is more than formal. Global analysis and differential geometry naturally give rise to function spaces equipped with an infinite-dimensional geometric structure. For example, if $\mathbf{E}\longrightarrow M$ is a vector bundle over a Riemannian manifold, the Sobolev spaces $W^{m,p}(M,\mathbf{E})$ are Banach spaces and serve as configuration spaces for numerous variational problems. Other important examples include the spaces of smooth sections $\Gamma(\mathbf{E})$, spaces of smooth maps $C^\infty(M,N)$, diffeomorphism groups $\operatorname{Diff}(M)$, loop manifolds $C^\infty(S^1,M)$, and spaces of immersions or embeddings.

This viewpoint is particularly useful in the variational study of partial differential equations. Solutions are often characterized as critical points of an energy functional defined on an infinite-dimensional function space. It is also common to restrict such a functional to subsets with a manifold structure, such as unit spheres, Nehari manifolds, submanifolds defined by integral constraints, and level sets of differentiable functionals. The geometry of these spaces then provides the appropriate language for formulating criticality conditions, applying Lagrange multipliers, studying gradient flows, and analyzing constrained minimization problems.

Certain spaces of geometric objects can likewise be understood as infinite-dimensional manifolds. They include the space of Riemannian metrics on a fixed manifold, the space of connections on a principal bundle, various spaces of curves or surfaces, and some quotients by symmetry groups. This language is central to Morse theory, gauge theory, symplectic geometry, the calculus of variations, and the study of geometric equations.

The chapter develops the foundations of Banach manifolds that will be used in nonlinear analysis and variational methods. It includes their tangent spaces, submanifolds defined by differentiable constraints, local forms of immersions and submersions, the constant rank theorem, and the preimage theorem. It also examines the assumptions under which bump functions and smooth partitions of unity exist. This last question reveals an essential difference from finite dimensions: paracompactness of the manifold alone does not suffice if the model Banach space lacks the necessary smooth functions. The resulting framework covers many situations modeled on Sobolev spaces and prepares the later study of variational problems and mapping spaces.

\begin{semblanzaHistorica}{From finite-dimensional manifolds to configuration spaces}
The idea of modeling a manifold on a function space made differential geometry applicable to problems in which the entire unknown is a function. Eells and Palais developed global tools for these spaces and showed that many finite-dimensional constructions survive under precise functional assumptions. The essential difference is that, in infinite dimensions, existence of complements, bump functions, or suitable metrics is no longer automatic. This chapter therefore keeps the assumptions needed to apply each theorem explicit~\cite{Eells1966,Palais1968}.
\end{semblanzaHistorica}

\section{Banach manifolds}
\label{sec:variedades-banach}

We begin by studying manifolds modeled on Banach spaces. They retain the differential calculus and local theorems needed to treat function spaces as geometric objects. In particular, Sobolev spaces of sections, constraints defined by differentiable functionals, and many level sets arising in variational methods have a Banach manifold structure under suitable assumptions.

Unless expressly stated otherwise, the model spaces in this chapter are real Banach spaces, and differentiability is understood over the real field.

\begin{definition}\label{def:variedad-banach}\index{Banach manifold}\index{manifold!Banach}\glsadd{variedad-banach}
Let $(M,\tau)$ be a Hausdorff topological space, and let $k\in\mathbb N\cup\{\infty\}$, with $k\geq1$. A \textbf{Banach chart} on $M$ is a pair $(U,\varphi)$, where $U\subseteq M$ is nonempty and open and $\varphi\colon U\longrightarrow V$ is a homeomorphism from $U$ onto an open subset $V$ of a Banach space $X$.

Two charts $(U_\alpha,\varphi_\alpha)$ and $(U_\beta,\varphi_\beta)$, modeled on Banach spaces $X_\alpha$ and $X_\beta$, respectively, are \textbf{$C^k$--compatible} if $U_\alpha\cap U_\beta=\varnothing$ or if the coordinate change
\[
\varphi_\beta\circ\varphi_\alpha^{-1}\colon
\varphi_\alpha(U_\alpha\cap U_\beta)\longrightarrow
\varphi_\beta(U_\alpha\cap U_\beta)
\]
is a diffeomorphism of class $C^k$ between open subsets of Banach spaces. Differentiability is understood in the Fréchet sense, as in Definition~\ref{def:funcion-Ck-banach}, and diffeomorphism has the meaning given in Definition~\ref{def:difeomorfismo-banach}. When $k=\infty$, both the coordinate change and its inverse are required to be of class $C^r$ for every $r\in\mathbb N$.

A \textbf{Banach atlas of class $C^k$} on $M$ is a family of charts $\{(U_\alpha,\varphi_\alpha)\}_{\alpha\in I}$ covering $M$ whose charts are mutually $C^k$--compatible. A \textbf{Banach manifold without boundary of class $C^k$} is a Hausdorff topological space equipped with a maximal Banach atlas of class $C^k$. A Banach manifold of class $C^\infty$ is called a \textbf{smooth Banach manifold}.
\end{definition}

\begin{remark}[Topological conventions]
\label{obs:convencion-banach-moore}
Moore \cite[Section 1.3]{Moore2017GlobalAnalysis} assumes second countability
and a separable model on each connected component. Our definition does not
include these two restrictions: it allows nonseparable models and imposes
paracompactness or metrizability separately when a global result needs them.
These are not contradictory definitions, but conventions of different scope.
We keep the broader convention and do not infer smooth partitions of unity
from the mere existence of an atlas.
\end{remark}

In this chapter, and whenever Definition~\ref{def:variedad-banach} is invoked, the term \emph{Banach manifold} means a manifold without boundary unless expressly stated otherwise. We adopt the same convention for finite-dimensional manifolds occurring as domains or codomains of the maps under study. This assumption is incorporated into the local models: chart images are open in Banach spaces.

Banach manifolds with boundary also exist in infinite dimensions. Their charts allow relatively open subsets of half-spaces such as $X\times[0,\infty)$, with coordinate changes differentiable up to the boundary. The product $X\times[0,\infty)$ itself is an example, with boundary $X\times\{0\}$; see \cite[pp.~766--767]{Eells1966}. The local forms and tangent construction below use the preceding convention.

The boundary of a configuration space differs from the boundary of the manifold on which its functions are defined. If a Riemannian manifold $M$ has boundary, the space $W_0^{1,p}(M)$, with $1<p<\infty$, is still a Banach space and, as a manifold modeled on itself, has no boundary. The Dirichlet condition here defines a vector space of admissible functions. Constraints given by regular values of functionals will be described through embedded submanifolds of these spaces.

In the preceding definition, we allow, a priori, charts modeled on different Banach spaces. However, on each connected component of the manifold, the model space is determined up to Banach space isomorphism.

\begin{lemma}\label{lem:modelo-variedad-banach-conexa}
Let $M$ be a Banach manifold of class $C^k$, with $k\geq1$. If $M$ is connected, all Banach spaces modeling its charts are linearly and topologically isomorphic.
\end{lemma}

\begin{proof}
If $M=\varnothing$, there are no charts to compare. Suppose $M\neq\varnothing$. First observe that if two charts in the differentiable structure have nonempty intersection, their model spaces are isomorphic as Banach spaces. Let $(U,\varphi)$ and $(V,\psi)$ be two charts such that $U\cap V\neq\varnothing$, modeled on the Banach spaces $X$ and $Y$, respectively. Fix $x\in U\cap V$ and write $\xi:=\varphi(x)$. Since the coordinate change $\psi\circ\varphi^{-1}$ is a diffeomorphism of class at least $C^1$, the operator
\[
d(\psi\circ\varphi^{-1})_{\xi}\colon X\longrightarrow Y
\]
is linear and continuous. Differentiating the identities between the inverse coordinate changes gives
\[
d(\varphi\circ\psi^{-1})_{\psi(x)}
\circ d(\psi\circ\varphi^{-1})_{\xi}
=
\operatorname{id}_X,
\qquad
d(\psi\circ\varphi^{-1})_{\xi}
\circ d(\varphi\circ\psi^{-1})_{\psi(x)}
=
\operatorname{id}_Y.
\]
Thus $d(\psi\circ\varphi^{-1})_{\xi}$ is a continuous linear isomorphism with continuous inverse.

Now take a fixed chart $(U_0,\varphi_0)$, modeled on a Banach space $X_0$, and define
\[
A:=
\left\{
x\in M\ \middle|\
\text{there exists a chart around $x$ modeled on a space isomorphic to $X_0$}
\right\}.
\]
The set $A$ is nonempty, since $U_0\subseteq A$, and open: if $x\in A$ and $(U,\varphi)$ is a chart around $x$ modeled on a space isomorphic to $X_0$, the same chart works around every point of $U$, so $U\subseteq A$.

The set $M\setminus A$ is also open. Let $y\in M\setminus A$ and let $(V,\psi)$ be a chart around $y$ modeled on a Banach space $Y$. If there existed $z\in V\cap A$, there would be a chart around $z$ modeled on a space isomorphic to $X_0$. This chart would intersect $(V,\psi)$ in an open set containing $z$, and the first part of the proof would show that $Y$ is isomorphic to $X_0$. The chart $(V,\psi)$ would then show that $y\in A$, a contradiction. Thus $V\subseteq M\setminus A$.

Hence $A$ is open and closed. Since $M$ is connected and $A\neq\varnothing$, we have $A=M$. Finally, if $(V,\psi)$ is any chart, choose $y\in V$. Since $y\in A$, there exists another chart around $y$ whose model space is isomorphic to $X_0$; comparing the two charts in a neighborhood of $y$, the first part of the proof shows that the model space of $(V,\psi)$ is also isomorphic to $X_0$. This completes the proof.
\end{proof}

This result justifies speaking, when $M$ is connected, of the model space of a Banach manifold up to Banach space isomorphism. For this reason, a Banach manifold is often defined from the outset as a manifold modeled on a fixed Banach space $X$.

Before introducing the tangent space, we specify the regularity of maps between Banach manifolds.

\begin{definition}\label{def:mapeo-Cr-variedades-banach}
\index{map!Ck between Banach manifolds@of class $C^k$ between Banach manifolds}
Let $M$ and $N$ be Banach manifolds of classes $C^k$ and $C^\ell$, respectively, and let $r\in\mathbb N\cup\{\infty\}$ be such that $1\leq r\leq k$ and $1\leq r\leq\ell$. A map $F\colon M\longrightarrow N$ is of \textbf{class $C^r$ at $p\in M$} if there exist a chart $(U,\varphi)$ of $M$ around $p$, a chart $(V,\psi)$ of $N$ around $F(p)$, and an open neighborhood $U_0\subseteq U$ of $p$ such that $F(U_0)\subseteq V$ and
\[
\psi\circ F\circ\varphi^{-1}\colon
\varphi(U_0)\longrightarrow\psi(V)
\]
is of class $C^r$. The map $F$ is of \textbf{class $C^r$} if it is of this class at every point of $M$. When $r=\infty$, we say that $F$ is \textbf{smooth}.
\end{definition}

This regularity does not depend on the chosen charts. Let $(\widetilde U,\widetilde\varphi)$ and $(\widetilde V,\widetilde\psi)$ be other charts around $p$ and $F(p)$, respectively. The representation $\psi\circ F\circ\varphi^{-1}$ is continuous and $V\cap\widetilde V$ is a neighborhood of $F(p)$; we may therefore restrict the domain to a neighborhood of $p$ whose image under $F$ is contained in $V\cap\widetilde V$. Restricting also to $U\cap\widetilde U$ gives
\[
\widetilde\psi\circ F\circ\widetilde\varphi^{-1}
=
(\widetilde\psi\circ\psi^{-1})
\circ(\psi\circ F\circ\varphi^{-1})
\circ(\varphi\circ\widetilde\varphi^{-1}).
\]
The maps on either side of the coordinate representation of $F$ are coordinate changes of class at least $C^r$. The chain rule shows that one coordinate representation of $F$ is of class $C^r$ if and only if every other one is.

\begin{remark}[The structure of lower regularity]
\label{obs:atlas-regularidad-menor-banach}
Let $\mathcal A^k$ be the maximal atlas of a manifold of class $C^k$, and fix $r\in\mathbb N\cup\{\infty\}$ with $1\leq r\leq k$. The coordinate changes of $\mathcal A^k$ are also of class $C^r$. This atlas therefore determines a structure of class $C^r$ on the same topological space; denote its maximal atlas by $\mathcal A^r$. Its charts are those that are $C^r$--compatible with the charts of $\mathcal A^k$. Any two such charts are compatible with each other: near each point of an overlap, insert a chart of $\mathcal A^k$ and write their coordinate change as a composition of two diffeomorphisms of class $C^r$.

When $r<k$, the atlas $\mathcal A^r$ may contain more charts than $\mathcal A^k$. A chart in $\mathcal A^r$ will be called a \emph{chart of class $C^r$ on $M$}. The class $C^r$ of a map can be checked in these charts, since the argument for coordinate independence uses only coordinate changes of class $C^r$. This distinction allows us to construct coordinates adapted to a map with the regularity that the map possesses. In the results on local forms of maps of class $C^r$, the constructed charts will belong to the underlying structures of class $C^r$.

Since $r\geq1$, all these coordinate changes are of class $C^1$. By the chain rule, two curves have the same velocity in a chart of $\mathcal A^k$ if and only if they do in a chart of $\mathcal A^r$. Thus passage to the structure of lower regularity preserves the tangent space to be defined next.
\end{remark}

To study locally Lipschitz vector fields and their flows,
we will also use regularity intermediate between $C^1$ and $C^2$.

\begin{definition}[Maps and manifolds of class $C^{1,1}$]
\label{def:variedad-banach-C11}
\index{manifold!Banach of class C11@Banach of class $C^{1,1}$}
\index{map!C11@of class $C^{1,1}$}
Let $X,Y$ be Banach spaces and let $U\subseteq X$ be open. A map
$f\colon U\longrightarrow Y$ is of \textbf{class $C^{1,1}$} if it is of class
$C^1$ and, for each $x_0\in U$, there exist $r,L>0$ such that
$B_X(x_0,r)\subseteq U$ and
\[
 \|Df(x)-Df(y)\|_{\mathcal L(X,Y)}
 \leq L\|x-y\|_X,
 \qquad x,y\in B_X(x_0,r).
\]
The condition on the derivative is local: a single constant
$L$ valid on all of $U$ is not required. A diffeomorphism is of class $C^{1,1}$ if both it
and its inverse have this regularity.

Two Banach charts are \textbf{$C^{1,1}$--compatible} if their intersection
is empty or the coordinate change between them is a diffeomorphism of
class $C^{1,1}$. An \textbf{atlas of class $C^{1,1}$} is a covering by
mutually $C^{1,1}$--compatible charts. A \textbf{Banach manifold of
class $C^{1,1}$} is a Hausdorff space equipped with a maximal atlas of
that class. Since its coordinate changes are $C^1$, this atlas determines
a Banach manifold structure of class $C^1$ on the same space.
\end{definition}

\begin{proposition}[Stability of $C^{1,1}$ regularity]
\label{prop:calculo-C11-banach}
Every map of class $C^2$ between open subsets of Banach spaces is of class
$C^{1,1}$. The composition of two maps of class $C^{1,1}$, whenever
defined, is of class $C^{1,1}$. Moreover, if a diffeomorphism of class $C^1$
is of class $C^{1,1}$, so is its inverse.
\end{proposition}

\begin{proof}
If $f$ is of class $C^2$, continuity of $D^2f$ allows its norm to be bounded
by a constant $L$ on a sufficiently small ball around
each point of its domain.
Theorem~\ref{teo:desigualdades-valor-medio-direccionales}, applied to the
map $Df$ on this convex ball, gives
$\|Df(x)-Df(y)\|\leq L\|x-y\|$. This proves the first assertion.

Now let $f\colon U\subseteq X\longrightarrow Y$ and
$g\colon V\subseteq Y\longrightarrow Z$ be of class $C^{1,1}$, with
$f(U)\subseteq V$. Fix $x_0\in U$. We can choose a convex ball
$B\subseteq U$ around $x_0$ whose image lies in a convex
ball $C\subseteq V$, and constants $M_f,M_g,L_f,L_g>0$, such that
$\|Df\|\leq M_f$ on $B$, $\|Dg\|\leq M_g$ on $C$, and the derivatives are
Lipschitz there with constants $L_f,L_g$, respectively. The mean value
inequality gives $\|f(x)-f(y)\|\leq M_f\|x-y\|$ for $x,y\in B$.
By the chain rule,
\[
 D(g\circ f)(x)-D(g\circ f)(y)
 =[Dg(f(x))-Dg(f(y))]Df(x)
 +Dg(f(y))[Df(x)-Df(y)].
\]
Consequently,
\[
 \|D(g\circ f)(x)-D(g\circ f)(y)\|
 \leq (L_gM_f^2+M_gL_f)\|x-y\|,
 \qquad x,y\in B.
\]
The composition is therefore of class $C^{1,1}$.

Finally, let $f\colon U\longrightarrow V$ be a diffeomorphism of class
$C^1$ and suppose $f$ is of class $C^{1,1}$. Write $h=f^{-1}$.
Around $y_0\in V$, choose a convex ball $C$ such that $h(C)$ is
contained in a ball on which $Df$ is $L$--Lipschitz. Since $Dh$ is continuous,
we may shrink $C$ so that $\|Dh(y)\|\leq M$ there. Then
$\|h(y)-h(z)\|\leq M\|y-z\|$ for $y,z\in C$. The identity
$A^{-1}-B^{-1}=A^{-1}(B-A)B^{-1}$, applied to
$A=Df(h(y))$ and $B=Df(h(z))$, gives
\[
 Dh(y)-Dh(z)
 =Dh(y)[Df(h(z))-Df(h(y))]Dh(z).
\]
Consequently,
\[
 \|Dh(y)-Dh(z)\|\leq LM^3\|y-z\|,
 \qquad y,z\in C.
\]
This proves that $h$ is of class $C^{1,1}$.
\end{proof}

Stability under composition allows a $C^{1,1}$ atlas to be completed to
a unique maximal atlas by adding all charts compatible with the
given ones. To compare two added charts around a common point, insert
a chart of the original atlas and apply the preceding proposition
to the two coordinate changes. Likewise, a map between
$C^{1,1}$ manifolds is of class $C^{1,1}$ when its coordinate
representations are; the same composition of chart changes shows
that it suffices to check this in one chart around each point.

Every manifold of class at least $C^2$ thus has an underlying $C^{1,1}$
structure obtained by completing its atlas in this category. When working
at $C^{1,1}$ regularity, we use charts of this chosen atlas,
rather than all charts of the underlying $C^1$ atlas: a merely $C^1$ coordinate change
need not have locally Lipschitz differential. The admissible curves
and the tangent space construction below use only
the first-order chain rule, so they also apply to
this category through its underlying $C^1$ structure.

If $U\subseteq M$ is open, we equip it with the Banach manifold structure induced by that of $M$. Denote by $C^r(U,N)$ the set of maps of class $C^r$ from $U$ to $N$ and write $C^r(U):=C^r(U,\mathbb R)$. In particular, if $I\subseteq\mathbb R$ is an open interval, a curve $\gamma\colon I\longrightarrow M$ is of class $C^r$ if for each $t_0\in I$ there exist a chart $(U,\varphi)$ around $\gamma(t_0)$, modeled on $X$, and an open interval $J\subseteq I$ containing $t_0$ such that $\gamma(J)\subseteq U$ and $\varphi\circ\gamma\restriction_J\colon J\longrightarrow X$ is of class $C^r$. By chart independence, it suffices to verify this condition in one chart around each point of the curve.

We now introduce the notion of tangent space. In finite dimensions there are several equivalent constructions, including one based on curves and another formulated through derivations on germs of smooth functions. On Banach manifolds, these two constructions need not agree under the canonical identification, so it is useful to distinguish them before fixing our convention.

\begin{definition}\label{def:curva-admisible-banach}
Let $M$ be a Banach manifold of class $C^k$, with $k\geq1$, and let $p\in M$. An \textbf{admissible curve at $p$} is a map $\gamma\colon (-\varepsilon,\varepsilon)\longrightarrow M$ of class $C^1$ such that $\gamma(0)=p$, where $\varepsilon>0$ may depend on $\gamma$.

Two admissible curves $\gamma_1$ and $\gamma_2$ are \textbf{equivalent at $p$}, and we write $\gamma_1\sim_p\gamma_2$, if there exists a chart $(U,\varphi)$ around $p$ such that, after restricting both curves to a sufficiently small common interval around $0$,
\[
\frac{d}{dt}\Biggr|_{t=0}(\varphi\circ\gamma_1)(t)
=
\frac{d}{dt}\Biggr|_{t=0}(\varphi\circ\gamma_2)(t).
\]
\end{definition}

\begin{proposition}\label{prop:equivalencia-curvas-banach}
Let $M$ be a Banach manifold of class $C^k$, with $k\geq1$, and let $p\in M$. The condition in Definition~\ref{def:curva-admisible-banach} is independent of the chosen chart and defines an equivalence relation on the set of admissible curves at $p$.
\end{proposition}

\begin{proof}
Let $(U,\varphi)$ and $(V,\psi)$ be two charts around $p$, and let $\gamma_1$ and $\gamma_2$ be two admissible curves at $p$. Since $\gamma_1(0)=\gamma_2(0)=p$ and the curves are continuous, there exists $\delta>0$ such that $\gamma_j((-\delta,\delta))\subseteq U\cap V$ for $j\in\{1,2\}$. Near $0$, we have
\[
\psi\circ\gamma_j
=
(\psi\circ\varphi^{-1})\circ(\varphi\circ\gamma_j),
\],
and the chain rule gives
\[
(\psi\circ\gamma_j)'(0)
=
d(\psi\circ\varphi^{-1})_{\varphi(p)}
\bigl((\varphi\circ\gamma_j)'(0)\bigr).
\]
The operator $A_{\psi\varphi}(p):=d(\psi\circ\varphi^{-1})_{\varphi(p)}$ is a continuous linear isomorphism, with inverse $d(\varphi\circ\psi^{-1})_{\psi(p)}$. Consequently,
\[
(\varphi\circ\gamma_1)'(0)=(\varphi\circ\gamma_2)'(0)
\quad\Longleftrightarrow\quad
(\psi\circ\gamma_1)'(0)=(\psi\circ\gamma_2)'(0).
\]
Thus two curves have the same velocity in one chart around $p$ if and only if they do in every other chart. In particular, in the preceding definition, ``there exists a chart'' can equivalently be replaced by ``for every chart.''

Reflexivity and symmetry are immediate. To check transitivity, suppose $\gamma_1\sim_p\gamma_2$ and $\gamma_2\sim_p\gamma_3$, and choose any chart $(U,\varphi)$ around $p$. By chart independence,
\[
(\varphi\circ\gamma_1)'(0)
=
(\varphi\circ\gamma_2)'(0)
=
(\varphi\circ\gamma_3)'(0),
\],
so $\gamma_1\sim_p\gamma_3$.
\end{proof}

\begin{definition}\label{def:tangente-cinematico-banach}\index{kinematic tangent vector}
Let $M$ be a Banach manifold of class $C^k$, with $k\geq1$, and let $p\in M$. The set of equivalence classes of admissible curves at $p$ is called the \textbf{kinematic tangent space} of $M$ at $p$ and is denoted by $T_p^{\operatorname{kin}}M$. The class of an admissible curve $\gamma$ is denoted by $[\gamma]$.
\end{definition}

\begin{proposition}\label{prop:coordenadas-tangente-cinematico-banach}
Let $M$ be a Banach manifold of class $C^k$, with $k\geq1$, let $p\in M$, and let $(U,\varphi)$ be a chart around $p$ modeled on a Banach space $X$. The map
\[
\Theta_{\varphi,p}\colon T_p^{\operatorname{kin}}M\longrightarrow X,
\qquad
\Theta_{\varphi,p}([\gamma])
:=
\frac{d}{dt}\Biggr|_{t=0}(\varphi\circ\gamma)(t),
\]
is a bijection. The vector operations transported from $X$ through $\Theta_{\varphi,p}$ do not depend on the chosen chart. Each chart also determines a complete norm on $T_p^{\operatorname{kin}}M$, and norms determined by different charts satisfy the explicit inequalities in the proof.
\end{proposition}

\begin{proof}
For every admissible curve $\gamma$ at $p$, continuity of $\gamma$ allows its domain to be restricted so that its image lies in $U$. Chart independence, proved in Proposition~\ref{prop:equivalencia-curvas-banach}, shows that $\Theta_{\varphi,p}$ does not depend on the representative of $[\gamma]$, so it is well defined. If $\Theta_{\varphi,p}([\gamma_1])=\Theta_{\varphi,p}([\gamma_2])$, then $\gamma_1\sim_p\gamma_2$ by the definition of the relation itself; thus $\Theta_{\varphi,p}$ is injective.

To prove surjectivity, let $v\in X$ and write $a:=\varphi(p)$. Since $\varphi(U)$ is open, there exists $r>0$ such that $B_X(a,r)\subseteq\varphi(U)$. If $\varepsilon:=\frac{r}{1+\|v\|_X}$, then $a+tv\in\varphi(U)$ for $|t|<\varepsilon$. The curve $\gamma_v(t):=\varphi^{-1}(a+tv)$ is admissible at $p$ and satisfies $\Theta_{\varphi,p}([\gamma_v])=v$.

For $\xi,\eta\in T_p^{\operatorname{kin}}M$ and $\lambda\in\mathbb R$, define
\[
\xi+\eta
:=
\Theta_{\varphi,p}^{-1}
\bigl(\Theta_{\varphi,p}(\xi)+\Theta_{\varphi,p}(\eta)\bigr),
\qquad
\lambda\xi
:=
\Theta_{\varphi,p}^{-1}
\bigl(\lambda\Theta_{\varphi,p}(\xi)\bigr).
\]
Let $(V,\psi)$ be another chart around $p$ modeled on a Banach space $Y$. The chain rule gives
\[
\Theta_{\psi,p}
=
A_{\psi\varphi}(p)\circ\Theta_{\varphi,p},
\qquad
A_{\psi\varphi}(p)
:=
d(\psi\circ\varphi^{-1})_{\varphi(p)}.
\]
Since $A_{\psi\varphi}(p)$ is linear, the operations obtained through $\Theta_{\psi,p}$ agree with those obtained through $\Theta_{\varphi,p}$. In particular, $T_p^{\operatorname{kin}}M$ has a well-defined vector space structure; its zero vector is the class of the constant curve $t\longmapsto p$.

Finally, set $\|\xi\|_{\varphi,p}:=\|\Theta_{\varphi,p}(\xi)\|_X$. This norm is complete because $\Theta_{\varphi,p}$ is a linear isometry onto $X$. Moreover,
\[
\|\xi\|_{\psi,p}
\leq
\|A_{\psi\varphi}(p)\|_{\mathcal L(X,Y)}
\|\xi\|_{\varphi,p},
\qquad
\|\xi\|_{\varphi,p}
\leq
\|A_{\psi\varphi}(p)^{-1}\|_{\mathcal L(Y,X)}
\|\xi\|_{\psi,p}.
\]
Thus the preceding two inequalities show that the norms associated with different charts determine the same topology. The resulting linear structure and Banach topology are independent of the chart, although in general there is no distinguished norm.
\end{proof}

There is also an algebraic construction of the tangent space. To formulate it, we need $M$ to be smooth.

\begin{definition}\label{def:germenes-suaves-variedad-banach}
\index{germ of a function}
Let $M$ be a smooth Banach manifold and let $p\in M$. A \textbf{germ at $p$ of a smooth real function} is an equivalence class of pairs $(U,f)$, where $U\subseteq M$ is an open neighborhood of $p$ and $f\colon U\longrightarrow\mathbb R$ is smooth. Two pairs $(U,f)$ and $(V,g)$ are equivalent if there exists an open neighborhood $W\subseteq U\cap V$ of $p$ such that $f\restriction_W=g\restriction_W$. The set of these germs is denoted by $C^\infty_p(M)$.
\end{definition}

The preceding relation is an equivalence relation: reflexivity and symmetry are immediate, and transitivity follows by intersecting two neighborhoods on which the corresponding representatives agree. Addition, multiplication, and scalar multiplication in $C^\infty_p(M)$ are defined by choosing representatives, restricting them to a common neighborhood of $p$, and performing the usual operations there. These definitions do not depend on the representatives, since equivalent representatives agree on a sufficiently small neighborhood of $p$. Thus $C^\infty_p(M)$ is a commutative unital real algebra. Evaluation at $p$ is well defined and is an algebra homomorphism. As usual, we use the same symbol $f$ for a germ and any of its representatives when no confusion can arise, and simply write $f(p)$.

\begin{definition}\label{def:tangente-operacional-banach}\index{operational tangent vector}
Let $M$ be a smooth Banach manifold and let $p\in M$. An \textbf{operational tangent vector at $p$} is an algebraic linear functional
\[
D\colon C^\infty_p(M)\longrightarrow\mathbb R
\]
satisfying the Leibniz rule
\[
D(fg)=D(f)g(p)+f(p)D(g)
\]
for any $f,g\in C^\infty_p(M)$. The space of all operational tangent vectors at $p$ is denoted by $T_p^{\operatorname{op}}M$.
\end{definition}

No continuity condition is imposed on $D$, since no topology has been chosen on the germ algebra $C^\infty_p(M)$. With addition and scalar multiplication defined pointwise, $T_p^{\operatorname{op}}M$ is a vector space. Every derivation $D\in T_p^{\operatorname{op}}M$ annihilates constants: $D(1)=D(1\cdot1)=2D(1)$ gives $D(1)=0$, and linearity gives $D(c)=0$ for every $c\in\mathbb R$.

\begin{proposition}\label{prop:inyeccion-tangente-operacional-banach}
Let $M$ be a smooth Banach manifold and let $p\in M$. Every kinematic tangent vector determines an operational tangent vector. More precisely, there exists a canonical injective linear map
\[
\mathcal J_p\colon T_p^{\operatorname{kin}}M\longrightarrow T_p^{\operatorname{op}}M,
\qquad
\mathcal J_p([\gamma])=D_{[\gamma]},
\],
where
\[
D_{[\gamma]}f
:=
\frac{d}{dt}\Biggr|_{t=0}f(\gamma(t)).
\]
\end{proposition}

\begin{proof}
The expression is well defined with respect to the germ. If two representatives $f$ and $g$ agree on a neighborhood $W$ of $p$, continuity of $\gamma$ allows its domain to be restricted so that $\gamma(t)\in W$ near $0$; then $f(\gamma(t))=g(\gamma(t))$ near $0$, and the two derivatives agree.

It is also independent of the curve representing $[\gamma]$. Fix a chart $(U,\varphi)$ around $p$ modeled on a Banach space $X$ and write $a:=\varphi(p)$. After restricting the representative of $f$ and the chart to a common neighborhood of $p$, the chain rule gives
\begin{equation}
\label{eq:derivacion-cinematica-coordenadas}
D_{[\gamma]}f
=
D(f\circ\varphi^{-1})(a)
\bigl[\Theta_{\varphi,p}([\gamma])\bigr].
\end{equation}
The right-hand side depends only on the germ of $f$ and the class $[\gamma]$.

Linearity of $D_{[\gamma]}$ in $f$ follows from linearity of the derivative, and the usual product rule shows that
\[
D_{[\gamma]}(fg)
=
D_{[\gamma]}(f)g(p)+f(p)D_{[\gamma]}(g).
\]
Thus $D_{[\gamma]}\in T_p^{\operatorname{op}}M$.

Formula \eqref{eq:derivacion-cinematica-coordenadas}, linearity of $D(f\circ\varphi^{-1})(a)$, and the definition of the vector operations in $T_p^{\operatorname{kin}}M$ imply that $\mathcal J_p$ is linear. To prove injectivity, let $\xi\in T_p^{\operatorname{kin}}M$ be such that $\mathcal J_p(\xi)=0$ and set $v:=\Theta_{\varphi,p}(\xi)$. For each $\ell\in X'$, consider the smooth function $f_\ell\colon U\longrightarrow\mathbb R$ defined by $f_\ell(x):=\ell(\varphi(x)-a)$. Then
\[
0
=
\mathcal J_p(\xi)(f_\ell)
=
\ell(v).
\]
Since the continuous dual $X'$ separates points of $X$, it follows that $v=0$. Injectivity of $\Theta_{\varphi,p}$ implies $\xi=0$, so $\mathcal J_p$ is injective.
\end{proof}

In finite dimensions, the identification in
Definition~\ref{def canonica espacio tangente derivaciones en p}, completed
in Proposition~\ref{ejer:nociones-fundamentales-variedad-suave-frontera-recordemos-denota-anillo},
shows that $\mathcal J_p$ is an isomorphism. In infinite dimensions it may fail to
be surjective, since the Leibniz rule alone does not force an
algebraic derivation to depend only on the first derivative.

\begin{proposition}[The operational tangent space may contain nonkinematic vectors]
\label{prop:tangente-operacional-l2}
Let $E:=\ell^2(\mathbb N;\mathbb R)$ be regarded as a smooth Banach manifold, and let $0\in E$. Then the canonical map
\[
\mathcal J_0\colon T_0^{\operatorname{kin}}E\longrightarrow T_0^{\operatorname{op}}E
\]
is not surjective. Consequently, under the canonical identification given by $\mathcal J_0$, the kinematic tangent space is a proper subspace of the operational tangent space.
\end{proposition}

\begin{proof}
Let $\mathcal L_s^2(E;\mathbb R)$ be the space of continuous symmetric bilinear forms on $E$, equipped with the norm in Definition~\ref{def:operadores-multilineales-continuos}, and denote by $\mathcal F_s(E)$ the subspace of symmetric forms of finite rank. Proposition~\ref{prop:formas-bilineales-simetricas-hilbert} identifies each $B\in\mathcal L_s^2(E;\mathbb R)$ isometrically with a continuous self-adjoint linear operator $A_B\colon E\longrightarrow E$ through $B(x,y)=\langle A_Bx,y\rangle_E$ and shows that $\overline{\mathcal F_s(E)}$ corresponds to the space of compact self-adjoint operators.

Consider the form $B_0(x,y):=\langle x,y\rangle_E$. The operator associated with $B_0$ is the identity on $E$, which is not compact because the standard orthonormal basis $(e_n)_{n\in\mathbb N}$ lies in the unit ball and satisfies $\|e_n-e_m\|_E=\sqrt2$ for $n\neq m$. Thus
\[
B_0\notin\overline{\mathcal F_s(E)}.
\]

Form the quotient space
\[
\mathcal Q
:=
\mathcal L_s^2(E;\mathbb R)/\overline{\mathcal F_s(E)}
\]
and let $\pi\colon \mathcal L_s^2(E;\mathbb R)\longrightarrow\mathcal Q$ be the canonical projection. By Proposition~\ref{prop:norma-cociente-banach}, $\mathcal Q$ is a normed space and
\[
\|\pi(B)\|_{\mathcal Q}
=
\operatorname{dist}_{\mathcal L_s^2(E;\mathbb R)}
\bigl(B,\overline{\mathcal F_s(E)}\bigr).
\]
In particular, $b_0:=\pi(B_0)$ is nonzero and $\|b_0\|_{\mathcal Q}>0$.

On the one-dimensional subspace $\operatorname{span}\{b_0\}$,
define
\[
\lambda_0\colon \operatorname{span}\{b_0\}\longrightarrow\mathbb R,
\qquad
\lambda_0(tb_0):=t.
\]
This function is well defined because $b_0\neq0$: if
$tb_0=sb_0$, then $(t-s)b_0=0$ and therefore $t=s$. Moreover,
$\lambda_0$ is linear. By absolute homogeneity of the norm, $\|tb_0\|_{\mathcal Q}
=
|t|\,\|b_0\|_{\mathcal Q}$, and consequently
\[
|\lambda_0(tb_0)|
=
|t|
=
\frac{\|tb_0\|_{\mathcal Q}}{\|b_0\|_{\mathcal Q}}.
\]
In particular, for every $t\neq0$,
\[
\frac{|\lambda_0(tb_0)|}{\|tb_0\|_{\mathcal Q}}
=
\frac{1}{\|b_0\|_{\mathcal Q}}.
\]
Since every nonzero element of $\operatorname{span}\{b_0\}$ has the
form $tb_0$ with $t\neq0$, we conclude that
\[
\|\lambda_0\|
=
\sup_{\substack{x\in\operatorname{span}\{b_0\}\\x\neq0}}
\frac{|\lambda_0(x)|}{\|x\|_{\mathcal Q}}
=
\frac{1}{\|b_0\|_{\mathcal Q}}.
\]
Thus $\lambda_0$ is continuous. Corollary~\ref{cor:extension-hahn-banach-preserva-norma} gives a continuous linear extension $\widetilde\lambda\colon \mathcal Q\longrightarrow\mathbb R$ of $\lambda_0$. Define
\[
\Lambda:=\widetilde\lambda\circ\pi
\colon \mathcal L_s^2(E;\mathbb R)\longrightarrow\mathbb R.
\]
Then $\Lambda$ is linear and continuous, annihilates $\overline{\mathcal F_s(E)}$, and satisfies $\Lambda(B_0)=1$.

For a germ $f\in C^\infty_0(E)$, define
\[
\mathscr D(f):=\Lambda(D^2f(0)).
\]
This definition is independent of the representative, since two representatives of the same germ agree on a neighborhood of $0$ and have the same derivatives there. Theorem~\ref{teo:simetria-derivadas-superiores} ensures that $D^2f(0)\in\mathcal L_s^2(E;\mathbb R)$. Linearity of the second derivative and of $\Lambda$ shows that $\mathscr D$ is linear.

Let $f,g\in C^\infty_0(E)$ and choose representatives defined on a common open neighborhood of $0$. By Proposition~\ref{prop:leibniz-segunda-derivada},
\[
D^2(fg)(0)
=
f(0)D^2g(0)+g(0)D^2f(0)
+
Df(0)\otimes Dg(0)+Dg(0)\otimes Df(0).
\]
Corollary~\ref{cor:forma-simetrizada-rango-finito-hilbert} shows that
\[
C_{f,g}
:=
Df(0)\otimes Dg(0)+Dg(0)\otimes Df(0)
\]
is a symmetric bilinear form of rank at most two. Thus $C_{f,g}\in\mathcal F_s(E)\subseteq\overline{\mathcal F_s(E)}$ and $\Lambda(C_{f,g})=0$. Applying $\Lambda$ to the preceding formula gives
\[
\mathscr D(fg)
=
f(0)\mathscr D(g)+g(0)\mathscr D(f).
\]
Hence $\mathscr D\in T_0^{\operatorname{op}}E$.

It remains to prove that $\mathscr D$ is not in the image of $\mathcal J_0$. Consider $q\colon E\longrightarrow\mathbb R$, $q(x):=\|x\|_E^2$. For $x,h\in E$, we have $q(x+h)-q(x)=2\langle x,h\rangle_E+\|h\|_E^2$, whence
\[
Dq(x)[h]=2\langle x,h\rangle_E,
\qquad
D^2q(x)[h,k]=2\langle h,k\rangle_E.
\]
The second derivative is constant and all higher derivatives vanish; in particular, $q$ is smooth, $Dq(0)=0$, and $D^2q(0)=2B_0$. Consequently,
\[
\mathscr D(q)=2\Lambda(B_0)=2.
\]
In contrast, if $[\gamma]\in T_0^{\operatorname{kin}}E$ and $v:=\gamma'(0)$, the chain rule gives
\[
\mathcal J_0([\gamma])(q)
=
\frac{d}{dt}\Biggr|_{t=0}q(\gamma(t))
=
Dq(0)[v]
=
0.
\]
Thus $\mathscr D$ does not agree with $\mathcal J_0([\gamma])$ for any admissible curve $\gamma$, and $\mathcal J_0$ is not surjective.
\end{proof}

Since $\mathcal J_0$ is injective, it canonically identifies $T_0^{\operatorname{kin}}E$ with the subspace $\mathcal J_0(T_0^{\operatorname{kin}}E)$ of $T_0^{\operatorname{op}}E$; the proposition shows that this inclusion is proper. This means precisely that there are operational derivations that do not arise from velocities of curves and hence that the two tangent space constructions do not agree canonically in infinite dimensions. Failure of surjectivity of $\mathcal J_0$ does not by itself exclude the possible existence of an abstract noncanonical linear isomorphism between the two vector spaces: to rule out that possibility as well, one would need to compare, for example, their algebraic dimensions. Such a stronger assertion is unnecessary here; the geometrically relevant difference is that the natural identification provided by curves does not exhaust the space of derivations.

From now on, unless expressly stated otherwise, ``tangent space'' means the kinematic tangent space.

\begin{definition}\label{def:espacio-tangente-banach}\index{tangent space of a Banach manifold}
Let $M$ be a Banach manifold of class $C^k$, with $k\geq1$, and let $p\in M$. The \textbf{tangent space} of $M$ at $p$ is
\[
T_pM:=T_p^{\operatorname{kin}}M.
\]
Its elements are equivalence classes of admissible curves. If $(U,\varphi)$ is a chart around $p$ modeled on a Banach space $X$, the associated bicontinuous linear isomorphism is
\[
\Theta_{\varphi,p}\colon T_pM\longrightarrow X,
\qquad
[\gamma]\longmapsto
\frac{d}{dt}\Biggr|_{t=0}(\varphi\circ\gamma)(t).
\]
\end{definition}

\begin{definition}\label{def:diferencial-mapeo-variedades-banach}
Let $F\colon M\longrightarrow N$ be a map of class $C^1$ between Banach manifolds and let $p\in M$. The \textbf{differential of $F$ at $p$} is the map
\[
dF_p\colon T_pM\longrightarrow T_{F(p)}N,
\qquad
dF_p[\gamma]:=[F\circ\gamma].
\]
\end{definition}

The preceding definition makes sense because $F\circ\gamma$ is a curve of class $C^1$ and satisfies $(F\circ\gamma)(0)=F(p)$. Moreover, it is independent of the representative of $[\gamma]$. Let $(U,\varphi)$ be a chart of $M$ around $p$ and $(V,\psi)$ a chart of $N$ around $F(p)$. After restricting the domains, for every admissible curve $\gamma$ at $p$ we have
\[
(\psi\circ F\circ\gamma)'(0)
=
d(\psi\circ F\circ\varphi^{-1})_{\varphi(p)}
\bigl((\varphi\circ\gamma)'(0)\bigr).
\]
If $\gamma_1$ and $\gamma_2$ represent the same vector in $T_pM$, their velocities in the chart $\varphi$ agree, and the preceding formula shows that $F\circ\gamma_1$ and $F\circ\gamma_2$ have the same velocity in the chart $\psi$. Thus $[F\circ\gamma_1]=[F\circ\gamma_2]$.

The same formula gives the coordinate representation
\[
\Theta_{\psi,F(p)}
\circ dF_p
\circ\Theta_{\varphi,p}^{-1}
=
d(\psi\circ F\circ\varphi^{-1})_{\varphi(p)}.
\]
In particular, $dF_p$ is a continuous linear operator. The chain rule for maps between open subsets of Banach spaces also implies that, if $G\colon N\longrightarrow P$ is of class $C^1$, then
\[
d(G\circ F)_p
=
dG_{F(p)}\circ dF_p,
\qquad
d(\operatorname{id}_M)_p
=
\operatorname{id}_{T_pM}.
\]

\section{Complementability and linear decompositions}
\label{sec:complementabilidad-geometria-banach}

The local forms of submanifolds and maps of regular rank
depend on a linear property that goes unnoticed in finite dimensions.
A closed subspace of a Banach space need not admit a
closed complement. Before studying immersions and submersions,
it is therefore useful to make precise the relation between topological complements,
continuous projections, and one-sided inverses.

\begin{definition}
\label{def:subespacio-complementado-banach}
\index{complemented subspace}
Let $E$ be a Banach space and let $Y\subseteq E$ be a closed subspace.
We say that $Y$ is \textbf{complemented} in $E$ if there exists a
closed subspace $Z\subseteq E$ such that
\[
E=Y\oplus Z.
\]
This means that every $x\in E$ can be written uniquely as $x=y+z$,
with $y\in Y$ and $z\in Z$. The subspace $Z$ is called a
\textbf{topological complement} of $Y$.
\end{definition}

The word ``topological'' refers to the fact that both summands are Banach
spaces with the induced norms and the projections determined by the
decomposition are continuous. This last property is not an
additional assumption, as the following characterization shows.

\begin{proposition}[Characterization by continuous projections]
\label{prop:caracterizacion-subespacio-complementado-banach}
\label{prop:proyecciones-descomposicion-complementada}
Let $E$ be a Banach space and $Y\subseteq E$ a closed subspace. The
following assertions are equivalent.
\begin{enumerate}[label=(\alph*)]
\item The subspace $Y$ is complemented in $E$.
\item There exists an operator $P\in\mathcal L(E)$ such that $P^2=P$ and
$\operatorname{Ran}(P)=Y$.
\end{enumerate}
If \textup{(b)} holds, then
\[
E=Y\oplus\ker(P).
\]
\end{proposition}

\begin{proof}
Suppose $E=Y\oplus Z$, with $Y$ and $Z$ closed. Equip
$Y\times Z$ with the norm $\|(y,z)\|:=\|y\|+\|z\|$. The operator
\[
S\colon Y\times Z\longrightarrow E,
\qquad
S(y,z):=y+z,
\]
is linear, continuous, and bijective. Since $Y\times Z$ and $E$ are Banach
spaces, Corollary~\ref{teo:inverso-acotado} ensures that $S^{-1}$ is
continuous. If $\operatorname{pr}_Y$ and $\operatorname{pr}_Z$ denote the
coordinate projections of $Y\times Z$, the projections determined by
the direct sum are
\[
P_Y=\operatorname{pr}_Y\circ S^{-1},
\qquad
P_Z=\operatorname{pr}_Z\circ S^{-1}.
\]
Thus both are continuous. In particular, $P_Y^2=P_Y$ and
$\operatorname{Ran}(P_Y)=Y$.

Conversely, let $P\in\mathcal L(E)$ be a projection with range $Y$.
The subspace $\ker(P)$ is closed. For each $x\in E$, we have
\[
x=Px+(x-Px),
\],
where $Px\in Y$ and $P(x-Px)=0$. Thus
$E=Y+\ker(P)$. If $x\in Y\cap\ker(P)$, then $x=Py$ for some
$y\in E$ and
\[
x=Py=P^2y=P(Py)=Px=0.
\]
The sum is direct and $E=Y\oplus\ker(P)$.
\end{proof}

\begin{definition}[One-sided inverses]
\label{def:inversas-laterales-banach}
Let $E$ and $F$ be Banach spaces and let
$T\in\mathcal L(E,F)$. An operator $L\in\mathcal L(F,E)$ is a
\textbf{left inverse} of $T$ if $L\circ T=I_E$, and an operator
$R\in\mathcal L(F,E)$ is a \textbf{right inverse} of $T$ if
$T\circ R=I_F$. Both are called \textbf{one-sided inverses} of $T$.
\end{definition}

\begin{proposition}[One-sided inverses and complementability]
\label{prop:inversas-laterales-complementabilidad-banach}
Let $E$ and $F$ be Banach spaces and let $T\in\mathcal L(E,F)$.
\begin{enumerate}[label=(\alph*)]
\item There exists $L\in\mathcal L(F,E)$ such that $L\circ T=I_E$ if and only if
$T$ is injective and $\operatorname{Ran}(T)$ is a complemented subspace of
$F$.
\item There exists $R\in\mathcal L(F,E)$ such that $T\circ R=I_F$ if and only if
$T$ is surjective and $\ker(T)$ is a complemented subspace of $E$.
\end{enumerate}
In the first case, $T\colon E\longrightarrow\operatorname{Ran}(T)$ is a
Banach space isomorphism. In the second, if
$E=\ker(T)\oplus E_1$, then
$T\restriction_{E_1}\colon E_1\longrightarrow F$ is an isomorphism.
\end{proposition}

\begin{proof}
First suppose $L\circ T=I_E$. Then $T$ is injective. The
operator $P:=T\circ L\in\mathcal L(F)$ satisfies
\[
P^2=T(LT)L=TL=P.
\]
Moreover, $P(F)\subseteq\operatorname{Ran}(T)$ and
$P(Tx)=T(LT)x=Tx$ for every $x\in E$. Consequently,
$\operatorname{Ran}(T)=\operatorname{Ran}(P)=\ker(I_F-P)$. This subspace
is closed, and
Proposition~\ref{prop:caracterizacion-subespacio-complementado-banach}
shows that this range is complemented.

Now suppose $T$ is injective and
$F=\operatorname{Ran}(T)\oplus Z$. Since the range is closed,
$T\colon E\longrightarrow\operatorname{Ran}(T)$ is a continuous bijective linear operator
between Banach spaces.
Corollary~\ref{teo:inverso-acotado} implies that its inverse
$T^{-1}\colon \operatorname{Ran}(T)\longrightarrow E$ is continuous. If
$P_{\operatorname{Ran}(T)}\colon F\longrightarrow\operatorname{Ran}(T)$ is the
continuous projection associated with the decomposition, then
\[
L:=T^{-1}\circ P_{\operatorname{Ran}(T)}
\]
belongs to $\mathcal L(F,E)$ and satisfies $L\circ T=I_E$. This proves
\textup{(a)}.

Suppose $T\circ R=I_F$. In particular, $T$ is surjective. The
operator $Q:=R\circ T\in\mathcal L(E)$ is a projection, since
\[
Q^2=R(TR)T=RT=Q.
\]
If $x\in\ker(Q)$, then $Tx=T(Qx)=0$. The reverse inclusion is
immediate, so $\ker(Q)=\ker(T)$. Likewise,
$\operatorname{Ran}(Q)=\operatorname{Ran}(R)$. Since $Q$ is a
continuous projection, its range and kernel are closed, and
Proposition~\ref{prop:caracterizacion-subespacio-complementado-banach}
gives the topological decomposition
\[
E=\ker(T)\oplus\operatorname{Ran}(R).
\]

Conversely, suppose $T$ is surjective and
$E=\ker(T)\oplus E_1$. The restriction
$T_1:=T\restriction_{E_1}\colon E_1\longrightarrow F$ is injective. It is also
surjective: given $y\in F$, choose $x\in E$ with $Tx=y$ and write
$x=k+x_1$, with $k\in\ker(T)$ and $x_1\in E_1$; then $T_1x_1=y$.
Corollary~\ref{teo:inverso-acotado} shows that
$T_1^{-1}\colon F\longrightarrow E_1$ is continuous. Composing this inverse
with the inclusion $E_1\hookrightarrow E$ gives a continuous right inverse
of $T$.
\end{proof}

One-sided inverses persist under sufficiently small perturbations.
This property explains why the local forms encountered below
persist in a neighborhood of the point under consideration.

\begin{corollary}[Stability of one-sided inverses]
\label{cor:estabilidad-inversas-laterales-banach}
Let $E$ and $F$ be Banach spaces and let $T,S\in\mathcal L(E,F)$.
\begin{enumerate}[label=(\alph*)]
\item If $L\in\mathcal L(F,E)$ satisfies $LT=I_E$ and
$\|L(S-T)\|_{\mathcal L(E)}<1$, then $S$ admits the continuous left inverse
\[
(LS)^{-1}L.
\]
\item If $R\in\mathcal L(F,E)$ satisfies $TR=I_F$ and
$\|(S-T)R\|_{\mathcal L(F)}<1$, then $S$ admits the continuous right inverse
\[
R(SR)^{-1}.
\]
\end{enumerate}
\end{corollary}

\begin{proof}
In the first case,
$LS=I_E+L(S-T)$ is invertible by
Theorem~\ref{teo:serie-de-neumann}; hence
$(LS)^{-1}LS=I_E$. In the second,
$SR=I_F+(S-T)R$ is invertible and
$SR(SR)^{-1}=I_F$.
\end{proof}

\begin{proposition}[Complements in the finite-dimensional and
Hilbert space cases]
\label{prop:complementabilidad-dimension-finita-hilbert}
Let $E$ be a Banach space.
\begin{enumerate}[label=(\alph*)]
\item Every finite-dimensional subspace of $E$ is complemented.
\item Every closed subspace of finite codimension in $E$ is
complemented.
\item If $E$ is a Hilbert space, every closed subspace of $E$ is
complemented.
\end{enumerate}
\end{proposition}

\begin{proof}
Let $Y\subseteq E$ be finite-dimensional and choose a basis
$y_1,\ldots,y_m$ of $Y$. The coordinate functionals
$\ell_j\colon Y\longrightarrow\mathbb R$ are continuous. By
Corollary~\ref{cor:extension-hahn-banach-preserva-norma}, each $\ell_j$ has
a continuous extension $\widetilde\ell_j\in E'$. The operator
\[
P(x):=\displaystyle\sum_{j=1}^{m}\widetilde\ell_j(x)y_j
\]
is a continuous projection onto $Y$.
Proposition~\ref{prop:caracterizacion-subespacio-complementado-banach}
proves \textup{(a)}.

Now let $Y\subseteq E$ be closed and of codimension $m<\infty$. Choose
$e_1,\ldots,e_m\in E$ whose classes form an algebraic basis of $E/Y$ and
set $Z:=\operatorname{span}\{e_1,\ldots,e_m\}$. Then $Z$ is closed,
$Y\cap Z=\{0\}$, and every element of $E$ belongs to $Y+Z$. Thus
$E=Y\oplus Z$.

Finally, if $E$ is a Hilbert space and $Y$ is closed,
Corollary~\ref{cor:descomposicion-ortogonal-hilbert} gives
$E=Y\oplus Y^\perp$.
\end{proof}

\section{Banach submanifolds}
\label{sec:subvariedades-banach}

We now turn to the notion of submanifold. In finite dimensions, a subset is a submanifold if it can locally be straightened by a change of coordinates. This principle remains valid in Banach spaces, but an additional topological condition arises: the subspace describing the tangent directions must be complemented. In infinite dimensions, a closed subspace need not have a closed complement, so this property must be part of the definition.

\begin{definition}\label{def:subvariedad-encajada-banach}\index{embedded submanifold}\index{manifold!Banach submanifold}\glsadd{subvariedad-banach}
Let $M$ be a Banach manifold of class $C^k$ and let $1\leq r\leq k$. A subset $N\subseteq M$ is an \textbf{embedded submanifold of class $C^r$} if, for each $p\in N$, there exist a chart $(U,\varphi)$ of class $C^r$ on $M$ around $p$, in the sense of Remark~\ref{obs:atlas-regularidad-menor-banach}, with $\varphi(p)=0$, and a topological decomposition $X=Y\oplus Z$ of the chart's model space into closed subspaces such that, upon identifying $X$ with $Y\times Z$ through $(y,z)\mapsto y+z$,
\[
\varphi(U\cap N)=\varphi(U)\cap(Y\times\{0\}).
\]
A chart of class $C^r$ with these properties is called an \textbf{adapted chart} for $N$ at $M$. The subspace $Y$ is the local model space of $N$ in this chart.
\end{definition}

The condition in Definition~\ref{def:subvariedad-encajada-banach} does more than describe the local position of $N$ inside $M$: it canonically equips $N$ with a Banach manifold structure.

\begin{proposition}[Induced manifold structure on a submanifold]
\label{prop:subvariedad-hereda-estructura-banach}
Let $N$ be an embedded submanifold of class $C^r$ of a Banach manifold $M$. Then $N$, with the subspace topology, is a Banach manifold of class $C^r$. If $(U,\varphi)$ is an adapted chart, $X=Y\oplus Z$ is the corresponding decomposition, and $\pi_Y\colon X\longrightarrow Y$ is the canonical projection, then
\[
\varphi_N:=\pi_Y\circ\varphi\restriction_{U\cap N}
\]
is a chart on $N$. The canonical inclusion $\iota\colon N\hookrightarrow M$ is of class $C^r$ and, for each $p\in N$, the differential $d\iota_p$ identifies $T_pN$ with a closed complemented subspace of $T_pM$.
\end{proposition}

\begin{proof}
Let $(U,\varphi)$ be an adapted chart around $p\in N$, with $\varphi(p)=0$ and $X=Y\oplus Z$. Define
\[
V_Y:=\{y\in Y\mid (y,0)\in\varphi(U)\}.
\]
The linear inclusion $\iota_Y\colon Y\longrightarrow X$, given by $\iota_Y(y)=(y,0)$, is continuous, and $V_Y=\iota_Y^{-1}(\varphi(U))$. Since $\varphi(U)$ is open in $X$, the set $V_Y$ is open in $Y$.

By the adapted chart condition, every point of $\varphi(U\cap N)$ has the form $(y,0)$ with $y\in V_Y$. Consequently, $\varphi_N\colon U\cap N\longrightarrow V_Y$ is bijective. Its inverse is given by
\[
\varphi_N^{-1}(y)=\varphi^{-1}(y,0),
\qquad y\in V_Y.
\]
Continuity of $\varphi_N$ follows from continuity of $\varphi\restriction_{U\cap N}$ and Proposition~\ref{prop:proyecciones-descomposicion-complementada}; continuity of $\varphi_N^{-1}$ follows from continuity of $\varphi^{-1}$ and $\iota_Y$. Thus $\varphi_N$ is a homeomorphism for the subspace topology on $N$.

Consider two adapted charts $(U_1,\varphi_1)$ and $(U_2,\varphi_2)$ with decompositions $X_1=Y_1\oplus Z_1$ and $X_2=Y_2\oplus Z_2$. On the corresponding overlap, the coordinate change between the induced charts on $N$ can be written as
\[
\varphi_{N,2}\circ\varphi_{N,1}^{-1}
=
\pi_{Y_2}\circ\varphi_2\circ\varphi_1^{-1}\circ\iota_{Y_1}.
\]
The inclusions $\iota_{Y_1}$ and projections $\pi_{Y_2}$ are continuous and linear, while $\varphi_2\circ\varphi_1^{-1}$ is of class $C^r$. The chain rule shows that the coordinate change is of class $C^r$. Interchanging the indices gives the same regularity for its inverse. Thus the charts $\varphi_N$ form a Banach atlas of class $C^r$ on $N$. Since $N$ is a subspace of a Hausdorff space, it is also Hausdorff, establishing its Banach manifold structure.

In the coordinates of an adapted chart, the canonical inclusion $\iota\colon N\hookrightarrow M$ is represented by the continuous linear operator
\[
\iota_Y\colon Y\longrightarrow Y\oplus Z,
\qquad
\iota_Y(y)=(y,0).
\]
Thus $\iota$ is of class $C^r$. Moreover, $d\iota_p$ is represented by $\iota_Y$, whose image is $Y\times\{0\}$. This subspace is closed and complemented by $\{0\}\times Z$, so $T_pN$ is identified with a closed complemented subspace of $T_pM$.
\end{proof}

Adapted charts also allow us to recognize when a map with image in $N$ is differentiable as a map into the submanifold itself.

\begin{lemma}[Maps into a submanifold]
\label{lem:corestriccion-subvariedad-banach}
Let $N\subseteq M$ be an embedded submanifold of class $C^r$, let $P$ be a Banach manifold of class at least $C^s$, and suppose $1\leq s\leq r$. If $F\colon P\longrightarrow M$ has image contained in $N$, then $F$ is of class $C^s$ if and only if its corestriction $\widehat F\colon P\longrightarrow N$ is of class $C^s$.
\end{lemma}

\begin{proof}
If $\widehat F$ is of class $C^s$, so is $F=\iota\circ\widehat F$, since the inclusion $\iota\colon N\longrightarrow M$ is of class $C^r$. Conversely, suppose $F$ is of class $C^s$. Since $N$ has the subspace topology, the corestriction $\widehat F$ is continuous. Given $p\in P$, choose an adapted chart $(U,\varphi)$ of $M$ around $F(p)$, with decomposition $X=Y\oplus Z$. On the open set $F^{-1}(U)$, the representation of $\widehat F$ in the induced chart $\varphi_N$ is
\[
\varphi_N\circ\widehat F=\pi_Y\circ\varphi\circ F.
\]
After introducing a chart of $P$, the right-hand side is a composition of maps of class $C^s$, since $\varphi$ is of class $C^r$ and $\pi_Y$ is continuous and linear. Consequently, $\widehat F$ is of class $C^s$ near $p$. Since $p$ was arbitrary, it is of this class on all of $P$.
\end{proof}

\begin{remark}\label{obs:subvariedad-grafica-banach}
By Theorem~\ref{teo:funcion-implicita-banach}, the description through adapted charts is equivalent to locally representing the subset, after a change of coordinates, as the graph of a map of class $C^r$ of the form $h\colon A\longrightarrow Z$, where $A\subseteq Y$ is open. The adapted chart description directly displays the tangent space and its complement in the ambient space.
\end{remark}

When the complement $Z$ can be chosen of finite dimension $m$, we say that $N$ has codimension $m$. In particular, a submanifold of codimension one is locally modeled on the kernel of a nonzero continuous linear functional. This situation arises naturally in regular level sets and variational constraints.

\section{Fundamental examples of Banach submanifolds}
\label{sec:ejemplos-subvariedades-banach}

The adapted charts in
Definition~\ref{def:subvariedad-encajada-banach} can be constructed explicitly
in several basic examples. These examples also show how
complementability determines the tangent geometry.

\begin{proposition}[Complemented subspaces and their translates]
\label{prop:subespacios-complementados-subvariedades-banach}
Let $E$ be a Banach space, let $Y\subseteq E$ be a closed complemented
subspace, and let $a\in E$. Then $a+Y$ is a smooth embedded
submanifold of $E$. If $\iota\colon a+Y\longrightarrow E$ is the inclusion and
$\Theta_{\operatorname{id}_E,x}\colon T_xE\longrightarrow E$ is the
bicontinuous linear isomorphism in
Proposition~\ref{prop:coordenadas-tangente-cinematico-banach}, then
\[
\Theta_{\operatorname{id}_E,x}
\bigl(d\iota_x(T_x(a+Y))\bigr)=Y
\]
for every $x\in a+Y$.
\end{proposition}

\begin{proof}
Choose a closed complement $Z$ such that $E=Y\oplus Z$. For each
$x\in a+Y$, the translation
\[
\varphi_x\colon E\longrightarrow E,
\qquad
\varphi_x(u):=u-x,
\]
is a smooth diffeomorphism, $\varphi_x(x)=0$, and
\[
\varphi_x(a+Y)=Y\times\{0\},
\],
where $E$ is identified with $Y\times Z$. Thus $\varphi_x$ is an adapted
chart. In these coordinates, the inclusion of $a+Y$ into $E$ is
represented by $y\mapsto(y,0)$, whose differential has image $Y$.
\end{proof}

\begin{proposition}[Tangent space to a graph]
\label{prop:graficas-subvariedades-banach}
Let $E$ and $F$ be Banach spaces, let $U\subseteq E$ be open, and let
$h\in C^r(U,F)$, with $r\geq1$. Then
\[
\operatorname{graf}(h):=\{(x,h(x))\in E\times F\mid x\in U\}
\]
is an embedded submanifold of class $C^r$ of the open set $U\times F$. Moreover,
if $\iota\colon \operatorname{graf}(h)\longrightarrow U\times F$ is the
inclusion, then
\[
\Theta_{\operatorname{id}_{U\times F},(x,h(x))}
\bigl(d\iota_{(x,h(x))}
(T_{(x,h(x))}\operatorname{graf}(h))\bigr)
=\{(v,Dh(x)v)\mid v\in E\}.
\]
\end{proposition}

\begin{proof}
The graph description of submanifolds in
Remark~\ref{obs:subvariedad-grafica-banach} can be seen directly.
Fix $x_0\in U$ and consider the diffeomorphism
\[
\Phi_{x_0}\colon U\times F\longrightarrow (U-x_0)\times F,
\qquad
\Phi_{x_0}(x,y):=(x-x_0,y-h(x)).
\]
The map $\Phi_{x_0}$ is of class $C^r$, its inverse is given by
\[
\Phi_{x_0}^{-1}(\xi,z)
=(\xi+x_0,z+h(\xi+x_0)),
\],
and it satisfies
\[
\Phi_{x_0}(x_0,h(x_0))=(0,0),
\qquad
\Phi_{x_0}(\operatorname{graf}(h))=(U-x_0)\times\{0\}.
\]
Thus $\Phi_{x_0}$ is a chart of class $C^r$ compatible with the underlying $C^r$ structure of $U\times F$ and adapted at $(x_0,h(x_0))$. The
differential of the parametrization $x\mapsto(x,h(x))$ at $x_0$ is given
by $v\mapsto(v,Dh(x_0)v)$, yielding the tangent space formula,
since $x_0$ was arbitrary. This space is complemented by
$\{0\}\times F$, since, for each $(u,w)\in E\times F$,
\[
(u,w)=(u,Dh(x_0)u)+(0,w-Dh(x_0)u).
\]
\end{proof}

\begin{example}[Hyperplanes and spheres]
\label{ej:hiperplanos-esferas-subvariedades-banach}
Let $E$ be a real Banach space.
\begin{enumerate}[label=(\alph*)]
\item If $\ell\in E'$ is nonzero and $c\in\mathbb R$, then
$\ell^{-1}(c)$ is a smooth embedded submanifold of codimension one.
Indeed, if $e\in E$ satisfies $\ell(e)=1$, then
$E=\ker(\ell)\oplus\mathbb Re$ and
\[
\ell^{-1}(c)=ce+\ker(\ell).
\]
Proposition~\ref{prop:subespacios-complementados-subvariedades-banach}
shows that this affine hyperplane is an embedded submanifold. If
$\iota_c\colon \ell^{-1}(c)\longrightarrow E$ is the inclusion, for every
$x\in\ell^{-1}(c)$ we have
\[
\Theta_{\operatorname{id}_E,x}
\bigl(d(\iota_c)_x(T_x\ell^{-1}(c))\bigr)=\ker(\ell).
\]

\item Let $r\geq1$ and suppose the norm on $E$ is of class $C^r$ on
$E\setminus\{0\}$. Set
\[
S_E:=\{x\in E\mid\|x\|_E=1\}.
\]
For each $x\in S_E$, homogeneity of the norm implies
\[
D\|\cdot\|_E(x)(x)
=
\frac{d}{dt}\Biggr|_{t=0}\|(1+t)x\|_E
=
\|x\|_E.
\]
Thus $D\|\cdot\|_E(x)\neq0$ on the unit sphere. To obtain
an adapted chart, fix $x\in S_E$ and set
$Y:=\ker(D\|\cdot\|_E(x))$. We have
$E=Y\oplus\mathbb Rx$. The map
\[
F_x\colon \mathcal O_x\longrightarrow\mathbb R,
\qquad
F_x(y,t):=\|x+y+tx\|_E-1,
\],
where $\mathcal O_x$ is a sufficiently small open neighborhood of
$(0,0)$ in $Y\times\mathbb R$ that $x+y+tx\neq0$, is of class $C^r$ and
satisfies $F_x(0,0)=0$ and
$D_2F_x(0,0)(s)=s$.
Theorem~\ref{teo:funcion-implicita-banach} locally represents the set
$F_x^{-1}(0)$ as the graph of a function of class $C^r$ on $Y$.
If $h$ denotes this function, the derivative formula for the implicit
function gives
\[
Dh(0)
=-D_2F_x(0,0)^{-1}\circ D_1F_x(0,0)=0,
\],
since $D_1F_x(0,0)$ is the restriction of
$D\|\cdot\|_E(x)$ to $Y$. In particular, the graph has tangent space
$Y\times\{0\}$ at the origin.
The coordinate change that subtracts this function in the second variable
straightens the graph onto $Y\times\{0\}$. Thus $S_E$ is an
embedded submanifold of class $C^r$. Differentiating the identity defining
the sphere gives, if
$\iota_S\colon S_E\longrightarrow E$ denotes the inclusion,
\[
\Theta_{\operatorname{id}_E,x}
\bigl(d(\iota_S)_x(T_xS_E)\bigr)
=\ker\bigl(D\|\cdot\|_E(x)\bigr).
\]
If $E$ is a real Hilbert space, we may take
$q(x):=\|x\|_E^2$; then $Dq(x)v=2\langle x,v\rangle_E$ and
$\Theta_{\operatorname{id}_E,x}
\bigl(d(\iota_S)_x(T_xS_E)\bigr)=x^\perp$.
\end{enumerate}
\end{example}

\begin{example}[The group of invertible operators]
\label{ej:grupo-lineal-general-banach}
Let $E$ be a Banach space. The set
\[
\operatorname{GL}(E):=\{A\in\mathcal L(E)\mid A\text{ is invertible}\}
\]
is open in the Banach space $\mathcal L(E)$ by
Proposition~\ref{prop:inversion-operadores-suave}. Thus
$\operatorname{GL}(E)$ is a smooth Banach manifold, and
\[
\Theta_{\operatorname{id}_{\operatorname{GL}(E)},A}
\colon T_A\operatorname{GL}(E)\longrightarrow\mathcal L(E)
\]
is a bicontinuous linear isomorphism for every
$A\in\operatorname{GL}(E)$. Operator composition is
continuous and bilinear and, by
Proposition~\ref{prop:diferenciabilidad-operadores-multilineales}, smooth.
The map $A\mapsto A^{-1}$ is also smooth and satisfies
\[
D(A\mapsto A^{-1})(A)H=-A^{-1}HA^{-1}.
\]
\end{example}

\section{Immersions, submersions, and embeddings}
\label{sec:inmersiones-submersiones-banach}

In finite dimensions, injectivity or surjectivity of the differential
alone provides the local inclusion and projection forms. On
Banach manifolds, one must also retain the topological information
contained in a one-sided inverse.
Proposition~\ref{prop:inversas-laterales-complementabilidad-banach} allows
this condition to be expressed through complemented kernels and images.
The resulting local forms are those commonly used in
infinite-dimensional differential geometry; more general
formulations may be found in
\cite{Gloeckner2015,Akin1978,Lang1999}.

\begin{definition}
\label{def:inmersion-submersion-split-banach}\glsadd{inmersion-banach}\glsadd{submersion-banach}
\index{immersion!between Banach manifolds}
\index{submersion!between Banach manifolds}
Let $M$ and $N$ be Banach manifolds of class $C^k$, let
$f\colon M\longrightarrow N$ be a map of class $C^r$, with $1\leq r\leq k$, and let
$p\in M$.
\begin{enumerate}[label=(\alph*)]
\item We say that $f$ is a \textbf{weak immersion at $p$} if
$df_p\colon T_pM\longrightarrow T_{f(p)}N$ is injective.

\item We say that $f$ is a \textbf{split immersion at $p$} if $df_p$
admits a continuous linear left inverse.

\item We say that $f$ is a \textbf{weak submersion at $p$} if $df_p$
is surjective.

\item We say that $f$ is a \textbf{split submersion at $p$} if $df_p$
admits a continuous linear right inverse.
\end{enumerate}
The map $f$ is a weak immersion, a split immersion, a weak
submersion, or a split submersion if it has the corresponding property at every
point of $M$. From now on, within the category of Banach manifolds,
\textbf{immersion} and \textbf{submersion} mean the
split versions unless expressly stated otherwise.
\end{definition}

The condition is intrinsic: differentials of coordinate changes
are continuous linear isomorphisms, and existence of a
one-sided inverse is preserved under composition on either side with isomorphisms.
Proposition~\ref{prop:inversas-laterales-complementabilidad-banach}
gives the characterizations
\[
\begin{split}
f\text{ is a split immersion at }p
&\quad\Longleftrightarrow\quad
df_p\text{ is injective and }\operatorname{Ran}(df_p)
\text{ is complemented},\\
f\text{ is a split submersion at }p
&\quad\Longleftrightarrow\quad
df_p\text{ is surjective and }\ker(df_p)
\text{ is complemented}.
\end{split}
\]
In the first equivalence, the image is automatically closed. Thus
mere injectivity, injectivity with closed image, and injectivity with
closed complemented image are distinct conditions.

\begin{theorem}[Local form of split immersions]
\label{teo:forma-local-inmersion-split-banach}
Let $f\colon M\longrightarrow N$ be a map of class $C^r$ between Banach
manifolds, with $r\geq1$, and let $p\in M$. The following assertions are
equivalent.
\begin{enumerate}[label=(\alph*)]
\item The map $f$ is a split immersion at $p$.
\item There exist charts of class $C^r$, $(U,\varphi)$ of $M$ around $p$ and $(V,\psi)$ of
$N$ around $f(p)$, Banach spaces $E$ and $Z$, and an identification
of the model space of $V$ with $E\oplus Z$ such that
$\varphi(p)=0$, $\psi(f(p))=(0,0)$, and
\[
(\psi\circ f\circ\varphi^{-1})(x)=(x,0)
\]
for every $x\in\varphi(U)$.
\end{enumerate}
In particular, if $f$ is a split immersion at $p$, after shrinking
$U$ the map $f\restriction_U$ is injective and is a split immersion at every
point of $U$.
\end{theorem}

\begin{proof}
Suppose \textup{(a)}. Take initial charts around $p$ and
$f(p)$. Composing them with translations reduces the problem to a map
$F\colon \Omega\subseteq E\longrightarrow\Theta\subseteq H$ of class $C^r$,
where $\Omega$ and $\Theta$ are open sets containing the origin, such that
$F(0)=0$ and
$A:=DF(0)$ admits a continuous left inverse. By
Proposition~\ref{prop:inversas-laterales-complementabilidad-banach}, the
subspace $R:=\operatorname{Ran}(A)$ is closed, and there exists a closed
subspace $Z\subseteq H$ such that $H=R\oplus Z$. Denote the
continuous projections by $P_R\colon H\longrightarrow R$ and $P_Z\colon H\longrightarrow Z$.

The operator $A\colon E\longrightarrow R$ is an isomorphism. Consider
\[
G\colon \Omega\longrightarrow R,
\qquad
G(x):=P_RF(x).
\]
We have $DG(0)=A$.
Theorem~\ref{teo:funcion-inversa-local-banach} gives open sets
$\Omega_0\subseteq\Omega$ and $W\subseteq R$ containing the origin such
that $G\restriction_{\Omega_0}\colon \Omega_0\longrightarrow W$ is a
diffeomorphism of class $C^r$. In the new coordinate $y=G(x)$, write
\[
F\bigl(G^{-1}(y)\bigr)=y+h(y),
\qquad
h(y):=P_ZF\bigl(G^{-1}(y)\bigr)\in Z.
\]

Since $h(0)=0$ and $\Theta$ is open, there exist open neighborhoods
$W_0\subseteq W$ and $Z_0\subseteq Z$ of the origin such that
\[
W_0\times Z_0\subseteq\Theta,
\qquad
h(W_0)\subseteq Z_0.
\]
Shrink $\Omega_0$ to $G^{-1}(W_0)$. The triangular map
\[
\Xi\colon W_0\times Z\longrightarrow W_0\times Z,
\qquad
\Xi(y,z):=(y,z-h(y)),
\]
is a diffeomorphism of class $C^r$, with inverse
$\Xi^{-1}(y,z)=(y,z+h(y))$. After restricting it to an open neighborhood
$W_0\times Z_0$ of the origin, its image $\Xi(W_0\times Z_0)$ is open and
$\Xi\restriction_{W_0\times Z_0}$ is a diffeomorphism onto that image.
Inclusion $h(W_0)\subseteq Z_0$ ensures that the points in the image of
$F$ under consideration lie in the domain of this new chart, giving
\[
\Xi\left(F\bigl(G^{-1}(y)\bigr)\right)=(y,0).
\]
Finally, the isomorphism $A\colon E\longrightarrow R$ allows us to postcompose the
coordinate $G$ with $A^{-1}$ and the first component of the codomain
coordinate with the same isomorphism. In this way, both first factors
are identified with the model space $E$, and the preceding expression becomes
$x\mapsto(x,0)$. Transporting these coordinate changes to $M$ and $N$
gives \textup{(b)}.

Conversely, in the coordinates of \textup{(b)}, the differential of $f$ at
$p$ is represented by the inclusion
\[
E\longrightarrow E\oplus Z,
\qquad
v\longmapsto(v,0),
\],
which admits projection onto $E$ as a left inverse. The chart differentials
are isomorphisms, so $df_p$ also admits a continuous left
inverse.

The local expression $(x,0)$ shows both injectivity of
$f\restriction_U$ and persistence of the split condition on $U$.
\end{proof}

\begin{theorem}[Local form of split submersions]
\label{teo:forma-local-submersion-split-banach}
Let $f\colon M\longrightarrow N$ be a map of class $C^r$ between Banach
manifolds, with $r\geq1$, and let $p\in M$. The following assertions are
equivalent.
\begin{enumerate}[label=(\alph*)]
\item The map $f$ is a split submersion at $p$.
\item There exist charts of class $C^r$, $(U,\varphi)$ of $M$ around $p$ and $(V,\psi)$ of
$N$ around $f(p)$, and a Banach space $K$, such that the model
space of $U$ is identified with $K\oplus F$, where $F$ is the model space
of $V$, and
\[
(\psi\circ f\circ\varphi^{-1})(k,y)=y
\]
for every $(k,y)\in\varphi(U)$.
\end{enumerate}
In particular, every split submersion is locally open, admits local
sections around every point of its image, and remains a split submersion
in a neighborhood of the point under consideration.
\end{theorem}

\begin{proof}
Suppose \textup{(a)} and work in centered charts, so the
map is represented by $F\colon \Omega\subseteq E\longrightarrow H$, with
$F(0)=0$. Let $T:=DF(0)$. By
Proposition~\ref{prop:inversas-laterales-complementabilidad-banach}, there exists
a closed subspace $E_1\subseteq E$ such that
\[
E=K\oplus E_1,
\qquad
K:=\ker(T),
\]
and $T_1:=T\restriction_{E_1}\colon E_1\longrightarrow H$ is an isomorphism.
Let $P_K\colon E\longrightarrow K$ be the continuous projection and define
\[
\Gamma\colon \Omega\longrightarrow K\times H,
\qquad
\Gamma(x):=(P_Kx,F(x)).
\]
If $x=k+e_1$, then
\[
D\Gamma(0)(k+e_1)=(k,T_1e_1).
\]
This operator is an isomorphism from $E$ onto $K\times H$, with inverse
$(k,y)\mapsto k+T_1^{-1}y$.
Theorem~\ref{teo:funcion-inversa-local-banach} shows that $\Gamma$ is a
diffeomorphism of class $C^r$ between neighborhoods of the origin. Since the second
component of $\Gamma$ is $F$, in the coordinate $(k,y)=\Gamma(x)$ we have
\[
F\bigl(\Gamma^{-1}(k,y)\bigr)=y.
\]
This proves \textup{(b)}.

If \textup{(b)} holds, the differential is represented by the
projection $K\oplus F\longrightarrow F$, which admits the right inverse
$y\mapsto(0,y)$. This proves the equivalence.

In these coordinates, the image of any sufficiently small open
product is open. Moreover, the map $y\mapsto(0,y)$, transported through
the charts, provides a local section. Persistence of the property follows
from the same normal form.
\end{proof}

\begin{definition}
\label{def:encaje-split-banach}
\index{embedding!between Banach manifolds}
Let $M$ and $N$ be Banach manifolds. A map
$f\colon M\longrightarrow N$ of class $C^r$, with $r\geq1$, is called a \textbf{split embedding} if it is a
split immersion, is injective, and the map
\[
f\colon M\longrightarrow f(M)
\]
is a homeomorphism when $f(M)$ is equipped with the subspace topology
from $N$.
\end{definition}

\begin{proposition}[Embeddings and embedded submanifolds]
\label{prop:encajes-subvariedades-banach}
Let $f\colon M\longrightarrow N$ be a map of class $C^r$ between Banach
manifolds, with $r\geq1$.
\begin{enumerate}[label=(\alph*)]
\item If $f$ is a split embedding, then $f(M)$ is an embedded submanifold
of class $C^r$ of $N$, and $f\colon M\longrightarrow f(M)$ is a diffeomorphism.
\item If $S\subseteq N$ is an embedded submanifold of class $C^r$, its
inclusion $S\hookrightarrow N$ is a split embedding.
\item If $f(M)$ is an embedded submanifold and
$f\colon M\longrightarrow f(M)$ is a diffeomorphism, then $f$ is a split
embedding.
\end{enumerate}
\end{proposition}

\begin{proof}
Suppose $f$ is a split embedding and fix $p\in M$. By
Theorem~\ref{teo:forma-local-inmersion-split-banach}, there exist neighborhoods
$U\subseteq M$ of $p$ and $V\subseteq N$ of $f(p)$, with charts
$\varphi\colon U\longrightarrow\varphi(U)\subseteq E$ and
$\psi\colon V\longrightarrow\psi(V)\subseteq E\oplus Z$, in which
\[
(\psi\circ f\circ\varphi^{-1})(x)=(x,0).
\]
Since $f\colon M\longrightarrow f(M)$ is a homeomorphism, $f(U)$ is open in
$f(M)$. Thus there exists an open set $O\subseteq N$ such that
$f(U)=O\cap f(M)$. The set $\psi(V\cap O)$ is an open subset of the space
$E\oplus Z$ containing the origin. Choose open neighborhoods
$W_E\subseteq E$ and $W_Z\subseteq Z$ of the origin so that
\[
W_E\subseteq\varphi(U),
\qquad
W_E\times W_Z\subseteq\psi(V\cap O).
\]
Define
\[
U_0:=\varphi^{-1}(W_E),
\qquad
V_0:=\psi^{-1}(W_E\times W_Z).
\]
The local form of $f$ shows that
$f(U_0)\subseteq V_0\cap f(M)$. Conversely, if
$q\in V_0\cap f(M)$, then $q\in O\cap f(M)=f(U)$. Writing
$q=f(x)$ with $x\in U$ gives
$\psi(q)=(\varphi(x),0)$. Since $\psi(q)\in W_E\times W_Z$, it follows that
$\varphi(x)\in W_E$ and therefore $x\in U_0$. We have proved the equality
$f(U_0)=V_0\cap f(M)$ and consequently
\[
\psi(V_0\cap f(M))
=W_E\times\{0\}
=\psi(V_0)\cap(E\times\{0\}).
\]
Thus $\psi\restriction_{V_0}$ is an adapted chart for $f(M)$. The same
local expression shows that the corestriction of $f$ and its inverse are of
class $C^r$.

The second assertion follows from
Proposition~\ref{prop:subvariedad-hereda-estructura-banach}: the inclusion is
of class $C^r$, is a homeomorphism onto its image, and its differential is the
inclusion of a closed complemented subspace. The third assertion follows
by composing the diffeomorphism $M\longrightarrow f(M)$ with the split
inclusion $f(M)\hookrightarrow N$.
\end{proof}

\begin{remark}[Phenomena specific to infinite dimensions]
\label{obs:fenomenos-inmersion-submersion-banach}
The conditions in
Definition~\ref{def:inmersion-submersion-split-banach} cannot be reduced to
injectivity or surjectivity of the differential.
\begin{enumerate}[label=(\alph*)]
\item The linear inclusion $I\colon \ell^1(\mathbb N)\longrightarrow\ell^2(\mathbb N)$ is continuous and
injective, but its image is not closed. If the first $N$ terms of $x^{(N)}$
are equal to $\frac{1}{N}$ and the remaining terms are zero, then
\[
\|x^{(N)}\|_{\ell^1(\mathbb N)}=1,
\qquad
\|Ix^{(N)}\|_{\ell^2(\mathbb N)}=\frac{1}{\sqrt N}\longrightarrow0.
\]
If the image were closed, Corollary~\ref{teo:inverso-acotado} applied to
$I\colon \ell^1(\mathbb N)\longrightarrow\operatorname{Ran}(I)$ would give a constant $C>0$
such that $\|x\|_{\ell^1(\mathbb N)}\leq C\|Ix\|_{\ell^2(\mathbb N)}$, contradicting this
sequence. Thus the image is not closed. In particular, $I$ admits no continuous left
inverse and, regarded as a
map between Banach manifolds, is a weak immersion that is not a
split immersion; indeed, by linearity, $dI_x=I$ at every
$x\in\ell^1(\mathbb N)$.

\item The subspace
$c_0:=\{x=(x_n)\in\ell^\infty(\mathbb N)\mid x_n\longrightarrow0\}$ is closed in
$\ell^\infty(\mathbb N)$, since a uniform limit of sequences tending to zero
also tends to zero. However, $c_0$ is not complemented; this
last assertion is the classical theorem of
Phillips. Consequently, the inclusion
$c_0\hookrightarrow\ell^\infty(\mathbb N)$ is a linear map that is a homeomorphism
onto its image and has injective differential with closed image, but is not
a split immersion. Its differential at every point is the inclusion itself.

Noncomplementability of $c_0$ in $\ell^\infty(\mathbb N)$ is
Phillips's Theorem~7.5~\cite[p.~539]{Phillips1940}. Its proof uses the
representation of the dual of $\ell^\infty(\mathbb N)$ by finitely additive measures
and a weak convergence argument that is not needed in the rest of the book.
Here it is used only to exhibit a closed subspace whose inclusion is not
split; no later result depends on that dual representation.

\item The canonical projection
\[
\pi\colon \ell^\infty(\mathbb N)\longrightarrow\ell^\infty(\mathbb N)/c_0
\]
is linear, continuous, and surjective. Since $c_0$ is closed, the quotient
$\ell^\infty(\mathbb N)/c_0$ is a Banach space by
Proposition~\ref{prop:norma-cociente-banach}. At each point, the differential
of the projection is $\pi$ and its kernel is $c_0$. If it admitted a continuous linear right
inverse,
Proposition~\ref{prop:inversas-laterales-complementabilidad-banach}
would imply that $\ker(\pi)=c_0$ is complemented in $\ell^\infty(\mathbb N)$, which
is false. Thus $\pi$ is a weak submersion that is not a split submersion.
\end{enumerate}
\end{remark}

\section{Constant rank and preimage theorems}
\label{sec:rango-constante-preimagen-banach}

For an operator between infinite-dimensional spaces, rank is not
adequately described by an integer. The local form corresponding to the
rank theorem requires simultaneous control of the kernel, the image, and
their complements. The following assumption expresses that, near the basepoint,
the images of the differentials retain a common complement in the
codomain.

\begin{theorem}[Constant rank with complemented decompositions]
\label{teo:rango-constante-split-abiertos-banach}
Let $E$ and $F$ be Banach spaces, let $\Omega\subseteq E$ be open, and let
$f\in C^r(\Omega,F)$, with $r\geq1$. Fix $a\in\Omega$ and set
\[
A:=Df(a),
\qquad
K:=\ker(A),
\qquad
R:=\operatorname{Ran}(A).
\]
Suppose there exist closed subspaces $X\subseteq E$ and $C\subseteq F$
such that
\[
E=K\oplus X,
\qquad
F=R\oplus C,
\]
and these direct sums are topological. Suppose, moreover, that there exists
an open neighborhood $\Omega_1\subseteq\Omega$ of $a$ for which
\begin{equation}
\label{eq:complemento-fijo-rango-constante-banach}
F=\operatorname{Ran}(Df(x))\oplus C
\end{equation}
is a topological direct sum for every $x\in\Omega_1$. Then there exist
open sets $W_R\subseteq R$ and $W_K\subseteq K$ containing the origin, an open set
$U_0\subseteq\Omega_1$ containing $a$, an open set $V_0\subseteq F$
containing $f(a)$, and diffeomorphisms of class $C^r$
\[
\alpha\colon U_0\longrightarrow W_R\times W_K,
\qquad
\beta\colon V_0\longrightarrow\beta(V_0)\subseteq R\times C,
\]
such that $\alpha(a)=(0,0)$, $\beta(f(a))=(0,0)$, and
\[
(\beta\circ f\circ\alpha^{-1})(y,k)=(y,0)
\]
for every $(y,k)\in W_R\times W_K$.
The open set $V_0$ can be chosen inside any previously specified open
neighborhood of $f(a)$.
\end{theorem}

\begin{proof}
Denote by $P_K\colon E\longrightarrow K$, $P_X\colon E\longrightarrow X$,
$P_R\colon F\longrightarrow R$, and $P_C\colon F\longrightarrow C$ the continuous
projections associated with the decompositions in the statement. Since
$K=\ker(A)$ and $R=\operatorname{Ran}(A)$, the restriction
\[
A_X:=A\restriction_X\colon X\longrightarrow R
\]
is linear and bijective. Corollary~\ref{teo:inverso-acotado} implies that
$A_X^{-1}\colon R\longrightarrow X$ is continuous.

Define
\[
G\colon \Omega_1\longrightarrow R\times K,
\qquad
G(x):=\bigl(P_R(f(x)-f(a)),P_K(x-a)\bigr).
\]
If $v=k+\xi$, with $k\in K$ and $\xi\in X$, then
\[
DG(a)v=(A_X\xi,k).
\]
This operator is an isomorphism from $E$ onto $R\times K$; its inverse is
$(y,k)\mapsto k+A_X^{-1}y$. By
Theorem~\ref{teo:funcion-inversa-local-banach}, after shrinking
$\Omega_1$ there exist an open set $U_1\subseteq\Omega_1$ containing $a$ and
an open set $W\subseteq R\times K$ containing $(0,0)$ such that
$G\restriction_{U_1}\colon U_1\longrightarrow W$ is a diffeomorphism of class
$C^r$. We may choose open balls $W_R\subseteq R$ and $W_K\subseteq K$,
centered at the origin, such that $W_R\times W_K\subseteq W$, and replace
$U_1$ by the preimage of this product.

In the coordinates $(y,k)=G(x)$, set
\[
\widetilde f(y,k)
:=
f\bigl(G^{-1}(y,k)\bigr)-f(a).
\]
The definition of $G$ implies $P_R\widetilde f(y,k)=y$. Thus
there exists a map $h\in C^r(W_R\times W_K,C)$ such that
\[
\widetilde f(y,k)=(y,h(y,k)).
\]
Fix $(y,k)\in W_R\times W_K$ and $\eta\in K$. We have
\[
D\widetilde f(y,k)(0,\eta)
=
(0,D_2h(y,k)\eta).
\]
This vector belongs to $C$. On the other hand, the chain rule and
invertibility of $DG^{-1}(y,k)$ show that
\[
\operatorname{Ran}(D\widetilde f(y,k))
=
\operatorname{Ran}\left(Df\bigl(G^{-1}(y,k)\bigr)\right).
\]
The direct sum in
\eqref{eq:complemento-fijo-rango-constante-banach} implies that the intersection
of this range with $C$ is trivial. It follows that
$D_2h(y,k)\eta=0$ for every $\eta\in K$, and consequently
\[
D_2h(y,k)=0.
\]

The ball $W_K$ can be chosen convex. For $y\in W_R$ and $k\in W_K$, the
integral mean value formula in
Theorem~\ref{teo:formula-integral-valor-medio} gives
\[
h(y,k)-h(y,0)
=
\int_0^1D_2h(y,tk)k\,dt
=0.
\]
Thus $h(y,k)=h(y,0)$. Define the triangular change of coordinates
\[
\Xi\colon W_R\times C\longrightarrow W_R\times C,
\qquad
\Xi(y,z):=(y,z-h(y,0)).
\]
Its inverse is given by $\Xi^{-1}(y,z)=(y,z+h(y,0))$, so $\Xi$
is a diffeomorphism of class $C^r$. Moreover,
\[
\Xi(\widetilde f(y,k))=(y,0).
\]
Since $h(0,0)=0$, choose an open ball $W_C\subseteq C$ centered at the
origin and shrink $W_R$, if necessary, so that
\[
h(y,0)\in W_C
\]
for every $y\in W_R$. Given an open neighborhood of $f(a)$ in
$F$ in advance, the balls $W_R$ and $W_C$ can also be chosen so that their image
under the sum isomorphism $R\times C\longrightarrow F$, translated by
$f(a)$, is contained in that neighborhood.

Take
\[
U_0:=G^{-1}(W_R\times W_K)
\]
and define $\alpha:=G\restriction_{U_0}$. Likewise, let
\[
V_0:=f(a)+\{y+z\mid y\in W_R,\ z\in W_C\}.
\]
This set is open in $F$ because the sum map
$R\times C\longrightarrow F$ is a topological isomorphism. The map
\[
\beta\colon V_0\longrightarrow\Xi(W_R\times W_C),
\qquad
\beta(v):=
\Xi\bigl(P_R(v-f(a)),P_C(v-f(a))\bigr),
\]
is a diffeomorphism of class $C^r$, satisfies
$\beta(f(a))=(0,0)$, and, for $(y,k)\in W_R\times W_K$, satisfies
\[
(\beta\circ f\circ\alpha^{-1})(y,k)=(y,0).
\]
These are the required diffeomorphisms.
\end{proof}

\begin{theorem}[Constant rank on Banach manifolds]
\label{teo:rango-constante-split-variedades-banach}
Let $M$ and $N$ be Banach manifolds and let $f\colon M\longrightarrow N$ be a map of
class $C^r$, with $r\geq1$. Fix $p\in M$. Suppose that, in
charts centered at $p$ and $f(p)$, a coordinate representation
\[
F\colon \Omega\subseteq E\longrightarrow H,
\qquad
F(0)=0,
\]
has the following property. If
\[
A:=DF(0),
\qquad
K:=\ker(A),
\qquad
R:=\operatorname{Ran}(A),
\],
there exist closed subspaces $X\subseteq E$ and $C\subseteq H$ such that
\[
E=K\oplus X,
\qquad
H=R\oplus C,
\]
as topological direct sums, and $C$ remains a topological
complement of $\operatorname{Ran}(DF(x))$ for every $x$ in a neighborhood of the
origin. Then there exist charts of class $C^r$ in which $f$ has the form
\[
(y,k)\longmapsto(y,0).
\]
In particular, locally $f$ factors as the projection
$R\times K\longrightarrow R$ followed by the inclusion
$R\hookrightarrow R\times C$. Its local image is an embedded submanifold, and
its local fibers are embedded submanifolds modeled on $K$.
\end{theorem}

\begin{proof}
Apply
Theorem~\ref{teo:rango-constante-split-abiertos-banach} to the coordinate
representation of $f$, choosing the open set $V_0$ inside the image of the initial
codomain chart. Composing the resulting diffeomorphisms with the
initial charts gives charts of class $C^r$, as in
Remark~\ref{obs:atlas-regularidad-menor-banach}. The assertions about the
image, fibers, and factorization follow directly from the expression
$(y,k)\mapsto(y,0)$.
\end{proof}

The common complement assumption in the preceding theorem replaces the
numerical constancy of rank used in finite dimensions. A particular fixed
complement need not be preserved under a nonlinear change of
coordinates. The conclusion, however, is intrinsic: existence
of charts in which the map has the form $(y,k)\mapsto(y,0)$ does not depend
on the initial charts used to construct them.

\begin{proposition}[Preimage of a submanifold under a submersion]
\label{prop:preimagen-subvariedad-submersion-banach}
Let $M$ and $N$ be Banach manifolds, let $f\colon M\longrightarrow N$ be a
split submersion of class $C^r$, with $r\geq1$, and let $S\subseteq N$ be an embedded
submanifold of class $C^r$. Then $f^{-1}(S)$ is an embedded submanifold of
class $C^r$ of $M$. The restriction
\[
f\restriction_{f^{-1}(S)}\colon f^{-1}(S)\longrightarrow S
\]
is a split submersion and, for each $p\in f^{-1}(S)$, if
$\iota\colon f^{-1}(S)\longrightarrow M$ and
$\jmath\colon S\longrightarrow N$ are the inclusions, then
\begin{equation}
\label{eq:tangente-preimagen-subvariedad-submersion-banach}
d\iota_p\bigl(T_p(f^{-1}(S))\bigr)
=
(df_p)^{-1}\bigl(d\jmath_{f(p)}(T_{f(p)}S)\bigr).
\end{equation}
\end{proposition}

\begin{proof}
If $f^{-1}(S)=\varnothing$, all assertions hold
vacuously. Thus suppose $f^{-1}(S)\neq\varnothing$.
Fix $p\in f^{-1}(S)$ and write $q:=f(p)$. Choose a chart
$(V,\psi)$ of $N$ adapted to $S$ at $q$. After translating the origin and
identifying the model space with a sum $Y\oplus Z$, we have
\[
\psi(V\cap S)=\psi(V)\cap(Y\times\{0\}).
\]
Take an initial chart of $M$ around $p$ and represent $f$ in
these coordinates, retaining the adapted chart $\psi$ in the codomain. The
differential of the representation admits a continuous right inverse. The
construction in the proof of
Theorem~\ref{teo:forma-local-submersion-split-banach}, which leaves the
codomain coordinate unchanged, gives a domain chart in which
\[
f(k,y,z)=(y,z),
\],
where the domain model space is identified with $K\oplus Y\oplus Z$.
In these coordinates,
\[
f^{-1}(S)=\{(k,y,0)\mid (k,y,0)\text{ belongs to the chart domain}\}.
\]
This set is the intersection of the domain with the complemented subspace
$K\oplus Y\oplus\{0\}$. Thus the chart is adapted to $f^{-1}(S)$.

The restriction of $f$ is represented by the projection
\[
K\oplus Y\longrightarrow Y,
\qquad
(k,y)\longmapsto y,
\],
which admits the right inverse $y\mapsto(0,y)$. Finally, the differential of
$f$ in these coordinates is $(k,y,z)\mapsto(y,z)$, and the preimage of
$Y\times\{0\}$ is $K\oplus Y\oplus\{0\}$. This proves
\eqref{eq:tangente-preimagen-subvariedad-submersion-banach}.
\end{proof}

\begin{definition}
\label{def:valor-regular-split-banach}\glsadd{valor-regular-banach}
\index{regular value!on Banach manifolds}
Let $M$ and $N$ be Banach manifolds and let $f\colon M\longrightarrow N$ be a map of
class $C^r$, with $r\geq1$. A point $q\in N$ is called a \textbf{split regular value}
of $f$ if, for each $p\in f^{-1}(q)$, the differential
\[
df_p\colon T_pM\longrightarrow T_qN
\]
admits a continuous linear right inverse. Equivalently, $df_p$ is
surjective and $\ker(df_p)$ is complemented in $T_pM$ for every
$p\in f^{-1}(q)$. If $f^{-1}(q)=\varnothing$, the condition is deemed
satisfied.
\end{definition}

\begin{theorem}[Preimage theorem]
\label{teo:preimagen-valor-regular-split-banach}
Let $M$ and $N$ be Banach manifolds, let $f\colon M\longrightarrow N$ be a map of
class $C^r$, with $r\geq1$, and let $q\in N$ be a split regular value of $f$.
Then $f^{-1}(q)$ is an embedded submanifold of class $C^r$ of $M$ and,
if $\iota_q\colon f^{-1}(q)\longrightarrow M$ is the inclusion,
\[
d(\iota_q)_p\bigl(T_p(f^{-1}(q))\bigr)=\ker(df_p)
\]
for every $p\in f^{-1}(q)$. If $N$ is finite-dimensional of dimension $m$,
then $f^{-1}(q)$ has codimension $m$.
\end{theorem}

\begin{proof}
If the fiber is empty, the local condition in
Definition~\ref{def:subvariedad-encajada-banach} holds vacuously;
we adopt the usual convention that the empty submanifold has the
indicated codimension. Suppose the fiber is nonempty and let
$p\in f^{-1}(q)$. By split regularity and
Theorem~\ref{teo:forma-local-submersion-split-banach}, there exist charts
centered at $p$ and $q$ in which $f$ has the form
\[
(k,y)\longmapsto y.
\]
Consequently, the fiber of $q$ corresponds in these coordinates to
$K\times\{0\}$, where $K$ represents $\ker(df_p)$. This subspace is
closed and complemented by $\{0\}\times F$, where $F$ is the model
space of the codomain chart. The isomorphism
$d\psi_q\colon T_qN\longrightarrow F$ transports the tangent directions
of the codomain to this second factor. Thus the resulting charts are
adapted charts for $f^{-1}(q)$. Its tangent space at $p$ is $K$, which under
the chart differentials corresponds precisely to $\ker(df_p)$.

If $N$ has dimension $m$, the normal complement
$\{0\}\times F$ has dimension $m$, since $d\psi_q$ is an isomorphism, and the
fiber has codimension $m$.
\end{proof}

\begin{corollary}[Regular values with finite-dimensional codomain]
\label{cor:valor-regular-codominio-finito-banach}
Let $M$ be a Banach manifold, let $N$ be a finite-dimensional manifold,
let $f\colon M\longrightarrow N$ be a map of class $C^r$, with $r\geq1$, and let
$q\in N$. Set
$m:=\dim(N)$. If
$df_p\colon T_pM\longrightarrow T_qN$ is surjective for every
$p\in f^{-1}(q)$, then $q$ is a split regular value. Consequently,
$f^{-1}(q)$ is an embedded submanifold of codimension $m$ and, if
$\iota_q\colon f^{-1}(q)\longrightarrow M$ is the inclusion, for every
$p\in f^{-1}(q)$,
\[
d(\iota_q)_p\bigl(T_p(f^{-1}(q))\bigr)=\ker(df_p).
\]
\end{corollary}

\begin{proof}
If $f^{-1}(q)=\varnothing$, the regularity assumption is vacuous and the
conclusion is immediate under the usual convention that the empty set is
an embedded submanifold of the indicated codimension; the assertion about
tangent spaces is also vacuous. Suppose $f^{-1}(q)\neq\varnothing$.
Fix $p\in f^{-1}(q)$ and choose a basis $e_1,\ldots,e_m$ of $T_qN$.
By surjectivity of $df_p$, there exist
$v_1,\ldots,v_m\in T_pM$ such that $df_p(v_j)=e_j$. The operator
\[
R_p\colon T_qN\longrightarrow T_pM,
\qquad
R_p\left(\displaystyle\sum_{j=1}^{m}a_je_j\right)
:=
\displaystyle\sum_{j=1}^{m}a_jv_j
\]
is linear and continuous because its domain is finite-dimensional, and satisfies
$df_p\circ R_p=I_{T_qN}$. Thus $q$ is a split regular value, and
Theorem~\ref{teo:preimagen-valor-regular-split-banach} gives the
conclusion.
\end{proof}

\begin{example}[Fibers of linear operators]
\label{ej:fibras-operadores-lineales-split-banach}
Let $E$ and $F$ be Banach spaces and let $T\in\mathcal L(E,F)$ be surjective.
If $\ker(T)$ is complemented, every $y\in F$ is a split regular value of
$T$. For any $x_0\in T^{-1}(y)$, we have
\[
T^{-1}(y)=x_0+\ker(T).
\]
If $\iota_y\colon T^{-1}(y)\longrightarrow E$ is the inclusion, then
\[
\Theta_{\operatorname{id}_E,x}
\bigl(d(\iota_y)_x(T_x(T^{-1}(y)))\bigr)=\ker(T)
\quad\text{for every }x\in T^{-1}(y).
\]
Indeed, for each $x\in T^{-1}(y)$, the translation
\[
\theta_x\colon \ker(T)\longrightarrow T^{-1}(y),
\qquad
\theta_x(k):=x+k,
\]
is an affine diffeomorphism, and its differential is the inclusion
$\ker(T)\hookrightarrow E$. Thus the tangent space of the fiber,
regarded as a Banach manifold modeled on $\ker(T)$, is identified
with $\ker(T)$. If the kernel is complemented, this inclusion admits a
continuous left inverse, and the fiber is an embedded submanifold by
Proposition~\ref{prop:subespacios-complementados-subvariedades-banach}.
If the kernel is not complemented, the fiber is still a closed affine
subspace and a Banach manifold in its own right. If the fiber were an
embedded submanifold of $E$, the differential of its inclusion would have image
exactly $\ker(T)$. Indeed, every differentiable curve $\gamma$ in
the fiber satisfies $T\circ\gamma=y$ and therefore
$T\gamma'(0)=0$; conversely, if $k\in\ker(T)$, the curve
$t\mapsto x+tk$ remains in the fiber and has velocity $k$ at $t=0$. Under the assumed embedded submanifold structure, this curve is also of class $C^1$ as a curve with values in the fiber, by Lemma~\ref{lem:corestriccion-subvariedad-banach}.
Proposition~\ref{prop:subvariedad-hereda-estructura-banach} would then force
this tangent image to be complemented in $E$, contradicting
the assumption. Consequently, the fiber is not an embedded
submanifold in the sense of
Definition~\ref{def:subvariedad-encajada-banach}.
\end{example}

\section{Regular level sets}

Level sets of real functionals are the most frequent specialization of Theorem~\ref{teo:preimagen-valor-regular-split-banach}. In this case, a nonvanishing differential automatically provides a right inverse and a one-dimensional complement to the kernel. We will develop this situation directly to obtain explicit adapted charts and fix the formulation to be used for variational constraints.
\label{sec:conjuntos-nivel-regulares-banach}

\begin{lemma}\label{lem:nucleos-funcionales-no-nulos-isomorfos}
Let $X$ be a real Banach space and let $f,g\in X'$ be nonzero functionals.
Then $\operatorname{ker}(f)$ and $\operatorname{ker}(g)$ are isomorphic as
Banach spaces.
\end{lemma}

\begin{proof}
First prove that there exists $a\in X$ such that $g(a)=1$ and $f(a)\neq0$. Since
$g\neq0$, there exists $b\in X$ such that $g(b)\neq0$. Replacing $b$ by $\frac{b}{g(b)}$,
we may assume $g(b)=1$. If $f(b)\neq0$, take $a=b$. Thus suppose
$f(b)=0$. Since $f\neq0$, there exists $x\in X$ such that $f(x)\neq0$.
Define $y:=x-g(x)b$. Then $g(y)=g(x)-g(x)g(b)=0$ and
$f(y)=f(x)-g(x)f(b)=f(x)\neq0$. Thus, taking $a:=b+y$, we have
$g(a)=g(b)+g(y)=1$ and $f(a)=f(b)+f(y)=f(y)\neq0$.

Now define $T\colon \operatorname{ker}(f)\longrightarrow \operatorname{ker}(g)$
by $T(h):=h-g(h)a$. First check that $T$ is well defined. If
$h\in\operatorname{ker}(f)$, then
$g(T(h))=g(h)-g(h)g(a)=g(h)-g(h)=0$. Thus
$T(h)\in\operatorname{ker}(g)$.

The operator $T$ is linear and continuous, since it is the restriction to
$\operatorname{ker}(f)$ of the continuous linear operator $x\mapsto x-g(x)a$
defined on $X$. We now prove that $T$ is bijective. To this end, define
$S\colon \operatorname{ker}(g)\longrightarrow \operatorname{ker}(f)$ by
$S(k):=k-\displaystyle\frac{f(k)}{f(a)}a$. This formula is well defined because
$f(a)\neq0$. Moreover, if $k\in\operatorname{ker}(g)$, then
$f(S(k))=f(k)-\displaystyle\frac{f(k)}{f(a)}f(a)=0$, so
$S(k)\in\operatorname{ker}(f)$. Thus $S$ is well defined. It is also linear and
continuous.

Let us check that $S$ and $T$ are inverses. Let $h\in\operatorname{ker}(f)$. Since
$T(h)=h-g(h)a$, we have
$f(T(h))=f(h)-g(h)f(a)=-g(h)f(a)$, because $f(h)=0$. Then

\[
S(T(h))
=
T(h)-\frac{f(T(h))}{f(a)}a
=
h-g(h)a-\frac{-g(h)f(a)}{f(a)}a
=
h.
\]

On the other hand, let $k\in\operatorname{ker}(g)$. Since
$S(k)=k-\displaystyle\frac{f(k)}{f(a)}a$, we have
$g(S(k))=g(k)-\displaystyle\frac{f(k)}{f(a)}g(a)=-\displaystyle\frac{f(k)}{f(a)}$, because $g(k)=0$ and
$g(a)=1$. Thus

\[
T(S(k))
=
S(k)-g(S(k))a
=
k-\frac{f(k)}{f(a)}a+\frac{f(k)}{f(a)}a
=
k.
\]

Thus $S=T^{-1}$. Consequently,
$T\colon \operatorname{ker}(f)\longrightarrow\operatorname{ker}(g)$ is a continuous linear bijection
with continuous linear inverse.

Finally, since $f$ and $g$ are continuous, $\operatorname{ker}(f)$ and
$\operatorname{ker}(g)$ are closed subspaces of $X$. Thus both are
Banach spaces with the norm induced from $X$. We conclude that
$\operatorname{ker}(f)$ and $\operatorname{ker}(g)$ are isomorphic as Banach
spaces.
\end{proof}

\begin{theorem}[Regular level sets in open subsets of Banach spaces]\label{teo:conjuntos-nivel-regulares-banach}\index{level set!regular}\index{manifold!regular level submanifold}\glsadd{nivel-regular}
Let $X$ be a real Banach space, let $\Omega\subseteq X$ be nonempty and open, and let $J\colon \Omega\longrightarrow\mathbb{R}$ be a function of class $C^{k}$, with $k\geq1$. Let $c\in\mathbb{R}$ and suppose the level set $V:=\{u\in\Omega\mid J(u)=c\}$ is regular, that is, $DJ(u)\neq0$ for every $u\in V$. For $u\in\Omega$, the identity chart determines the operator
\[
\mathcal I_u:=\Theta_{\operatorname{id}_{\Omega},u}
\colon T_u\Omega\longrightarrow X.
\]
Proposition~\ref{prop:coordenadas-tangente-cinematico-banach} proves that
$\mathcal I_u$ is linear, bijective, and bicontinuous. Then $V$ is an
embedded Banach submanifold of class $C^{k}$ and codimension $1$ of
$\Omega$. If $\iota\colon V\longrightarrow\Omega$ is the
inclusion, then
\[
\mathcal I_u\bigl(d\iota_u(T_uV)\bigr)=\operatorname{ker}(DJ(u))
\]
for every $u\in V$.
\end{theorem}

\begin{proof}
If $V=\varnothing$, the assertion is vacuous. Suppose $V\neq\varnothing$ and fix $u\in V$. Since $DJ(u)\neq0$, there exists $e\in X$ such that $DJ(u)(e)\neq0$. Let $Y_u:=\operatorname{ker}(DJ(u))$. Then $Y_u$ is closed and
\[
X=Y_u\oplus\mathbb Re.
\]
For each $x\in X$, we have
\[
x=
\left(x-\frac{DJ(u)(x)}{DJ(u)(e)}e\right)
+
\frac{DJ(u)(x)}{DJ(u)(e)}e,
\],
where the first summand belongs to $Y_u$ and the intersection $Y_u\cap\mathbb Re$ is trivial.

Consider
\[
\Phi(y,t):=J(u+y+te)-c
\]
on the open subset of $Y_u\times\mathbb R$ consisting of pairs $(y,t)$ such that $u+y+te\in\Omega$. Moreover, $\Phi(0,0)=0$ and $D_2\Phi(0,0)(1)=DJ(u)(e)\neq0$. By Theorem~\ref{teo:funcion-implicita-banach}, there exist open sets $U_0\subseteq Y_u$ and $I_0\subseteq\mathbb R$ containing $0$ and a function $h\colon U_0\longrightarrow I_0$ of class $C^k$, with $h(0)=0$, such that
\[
\Phi(y,t)=0
\quad\Longleftrightarrow\quad
t=h(y)
\]
for every $(y,t)\in U_0\times I_0$.

Let $W:=\{u+y+te\mid y\in U_0,\ t\in I_0\}$. The linear operator
$L_u\colon Y_u\times\mathbb R\longrightarrow X$,
$L_u(y,t):=y+te$, is a topological isomorphism by
Proposition~\ref{prop:proyecciones-descomposicion-complementada}; hence
$W$ is open in $X$. Define
\[
\varphi\colon W\longrightarrow Y_u\times\mathbb R,
\qquad
\varphi(u+y+te):=(y,t-h(y)).
\]
Its image is the open set
\[
\varphi(W)=\{(y,s)\in U_0\times\mathbb R\mid s+h(y)\in I_0\},
\],
and its inverse is given by $\varphi^{-1}(y,s)=u+y+(s+h(y))e$. Consequently, $\varphi$ is a diffeomorphism of class $C^k$, $\varphi(u)=(0,0)$, and
\[
\varphi(W\cap V)=\varphi(W)\cap(Y_u\times\{0\}).
\]
Thus $V$ is an embedded submanifold of class $C^k$ and codimension $1$.

To identify the tangent space in the natural coordinates of the open set $\Omega\subseteq X$, differentiate the identity $\Phi(y,h(y))=0$ at $y=0$. For each $\eta\in Y_u$, this gives
\[
0
=
D\Phi(0,0)\bigl(\eta,Dh(0)(\eta)\bigr)
=
DJ(u)\bigl(\eta+Dh(0)(\eta)e\bigr)
=
Dh(0)(\eta)DJ(u)(e).
\]
Since $DJ(u)(e)\neq0$, it follows that $Dh(0)=0$. The differential of the inverse adapted chart therefore sends $(\eta,0)$ to $\eta$. Proposition~\ref{prop:subvariedad-hereda-estructura-banach} thus gives
\[
\mathcal I_u\bigl(d\iota_u(T_uV)\bigr)
=Y_u=\operatorname{ker}(DJ(u)).
\]
Finally, if $u_0\in V$ is fixed, Lemma~\ref{lem:nucleos-funcionales-no-nulos-isomorfos} shows that $\operatorname{ker}(DJ(u))$ is isomorphic to $\operatorname{ker}(DJ(u_0))$ for every $u\in V$. Thus, even if $V$ is disconnected, all its local model spaces are isomorphic.
\end{proof}

\begin{lemma}[Differential of a restriction]\label{lem:diferencial-restriccion-subvariedad}\index{differential of a restriction}
Let $X$ be a Banach space, let $\Omega\subseteq X$ be open, let $V\subseteq\Omega$ be a Banach submanifold, and let $\varphi\in C^{1}(\Omega,\mathbb{R})$. Then $\varphi\restriction_{V}\in C^{1}(V,\mathbb{R})$ and, for each $u\in V$,
\[
d(\varphi\restriction_V)_u(v)
=D\varphi(u)\bigl(\mathcal I_u(d\iota_u(v))\bigr)
\]
for every $v\in T_uV$. In particular,
\[
d(\varphi\restriction_V)_u
=D\varphi(u)\circ\mathcal I_u\circ d\iota_u.
\]
\end{lemma}

\begin{proof}
Let $\iota\colon V\longrightarrow\Omega$ be the canonical inclusion. Then $\varphi\restriction_V=\varphi\circ\iota$. Proposition~\ref{prop:subvariedad-hereda-estructura-banach} shows that $\iota$ is of class $C^1$, and the chain rule gives
\[
d(\varphi\restriction_V)_u=d\varphi_u\circ d\iota_u.
\]
Under the bicontinuous isomorphism $\mathcal I_u\colon T_u\Omega\longrightarrow X$
defined in Theorem~\ref{teo:conjuntos-nivel-regulares-banach}, the
differential $d\iota_u$ is represented by the inclusion of
$\mathcal I_u(d\iota_u(T_uV))$ into $X$.
Moreover, in the identity chart of $\Omega$ and the usual coordinate on
$\mathbb R$, the differential $d\varphi_u$ is represented by the Fréchet
derivative $D\varphi(u)$; that is, $d\varphi_u=D\varphi(u)\circ\mathcal I_u$,
identifying $T_{\varphi(u)}\mathbb R$ with $\mathbb R$. Thus
\[
d(\varphi\restriction_V)_u(v)
=D\varphi(u)\bigl(\mathcal I_u(d\iota_u(v))\bigr)
\]
for every $v\in T_uV$, as required.
\end{proof}

\subsection{Positive functions with an integral constraint}

\begin{example}[Positive functions of integral one]
\label{ej:densidades-positivas-banach}
Let $(M,\mathbf g)$ be a nonempty compact Riemannian manifold without
boundary, and let $k\in\mathbb N_0$. Consider
\[
\mathcal D^k(M):=
\left\{\rho\in C^k(M,\mathbb R)\ \middle|\
\rho(x)>0\text{ for every }x\in M,\quad
\int_M\rho\,d\lambda_{\mathbf g}=1\right\}.
\]
This set is a smooth embedded submanifold of codimension one of
$C^k(M,\mathbb R)$. Its tangent space, identified by the inclusion
into the function space, is
\[
T_\rho\mathcal D^k(M)
\simeq
\left\{h\in C^k(M,\mathbb R)\ \middle|\
\int_Mh\,d\lambda_{\mathbf g}=0\right\}.
\]
Moreover, it admits both a global affine chart and a smooth global
parametrization by functions of integral zero through the normalized
exponential.
\end{example}

\begin{proof}
Use the algebra $\mathcal A=C^k(M,\mathbb R)$ with the norm of
Lemma~\ref{lem:algebra-ck-compacta}. Set
\[
V:=\lambda_{\mathbf g}(M),\qquad
\ell(f):=\int_Mf\,d\lambda_{\mathbf g},\qquad
\rho_*:=\frac{1}{V},\qquad Y:=\ker\ell.
\]
Compactness and nonemptiness give $0<V<\infty$. Moreover,
$|\ell(f)|\leq V\|f\|_\infty\leq V\|f\|_{\mathcal A}$, so
$\ell$ is continuous and linear. Since $\ell(\rho_*)=1$, the
operator $P(f):=f-\ell(f)\rho_*$ is a continuous projection onto
$Y$, and
\[
\mathcal A=Y\oplus\mathbb R\rho_*.
\]
In particular, $Y$ is closed and is a Banach space.

The set $\mathcal A_+$ of strictly positive functions is open by
Example~\ref{ej:funciones-elementales-algebras}. The translation
$\rho\mapsto\rho-\rho_*$ identifies $\mathcal D^k(M)$ with the
open subset
\[
\{y\in Y\mid\rho_*+y>0\}
\]
of $Y$. This is the global affine chart. To verify the embedded
submanifold property as well, fix $\rho\in\mathcal D^k(M)$.
The coordinate change on the ambient space
\[
\begin{split}
\Psi_\rho\colon\mathcal A&\longrightarrow Y\times\mathbb R,\\
\Psi_\rho(f)&:=
\bigl(f-\rho-(\ell(f)-1)\rho_*,\ \ell(f)-1\bigr)
\end{split}
\]
is a continuous affine isomorphism with inverse
\[
\Psi_\rho^{-1}(y,s)=\rho+y+s\rho_*.
\]
It sends $\rho$ to $(0,0)$ and, upon restriction to the open set
$\mathcal A_+$, sends $\mathcal D^k(M)$ exactly to the intersection
of its image with $Y\times\{0\}$. It is therefore an adapted chart.
Its inverse restricted to $s=0$ has differential $h\mapsto h$.
This proves the identification of the tangent space. It can also be
seen directly: if $h\in Y$ and $m_\rho:=\min_M\rho>0$, the curve
$\rho+th$ remains in $\mathcal D^k(M)$ when
\[
|t|<\frac{m_\rho}{2(1+\|h\|_\infty)}.
\]
Its velocity at zero is $h$. Conversely, differentiating the integral
constraint of any admissible curve gives $\ell(h)=0$.

We now construct nonlinear global coordinates. For $f\in Y$, define
\[
Z(f):=\int_Me^f\,d\lambda_{\mathbf g},\qquad
\Phi(f):=\frac{e^f}{Z(f)}.
\]
The function $e^f$ is positive and $Z(f)>0$. Integration is a
continuous functional and the exponential is smooth on $\mathcal A$;
therefore $Z$ and $\Phi$ are smooth. Moreover,
$\Phi(f)\in\mathcal D^k(M)$. For $\rho\in\mathcal D^k(M)$, set
\[
\Lambda(\rho):=\log\rho
-\frac{1}{V}\left(\int_M\log\rho\,d\lambda_{\mathbf g}\right)1.
\]
This map takes values in $Y$ and is smooth because the logarithm is
smooth on $\mathcal A_+$. If
$\displaystyle c_\rho:=\displaystyle\frac{1}{V}\int_M\log\rho\,d\lambda_{\mathbf g}$,
then $e^{\Lambda(\rho)}=e^{-c_\rho}\rho$ and
$Z(\Lambda(\rho))=e^{-c_\rho}$, so
$\Phi(\Lambda(\rho))=\rho$. On the other hand,
$\log\Phi(f)=f-(\log Z(f))1$. Since $\ell(f)=0$, subtracting
the average of this expression gives $\Lambda(\Phi(f))=f$.
Thus $\Phi$ and $\Lambda$ are smooth inverses.

The derivatives describe the effect of normalization. For $f,h\in Y$,
\[
\begin{split}
DZ(f)h&=\int_Me^fh\,d\lambda_{\mathbf g},\\
D\Phi(f)h&=
\Phi(f)\left(h-\int_M\Phi(f)h\,d\lambda_{\mathbf g}\right).
\end{split}
\]
The integral of the second expression is zero. If
$\rho\in\mathcal D^k(M)$ and $v\in Y$ is a tangent direction, then
\[
D\Lambda(\rho)v
=\frac{v}{\rho}
-\frac{1}{V}\left(\int_M\frac{v}{\rho}\,d\lambda_{\mathbf g}\right)1.
\]
These derivatives are inverse on the tangent spaces, either by the
chain rule applied to the preceding identities or by substitution into
the formulas. The integral constraint removes precisely the constant
direction that exponential normalization does not distinguish.
\end{proof}

\subsection{A manifold of conformal metrics}

Here a metric of class $C^k$ is a positive definite symmetric tensor
of class $C^k$. Regularity of its coefficients and smoothness of maps
between the Banach spaces are distinct properties.

\begin{example}[Metrics conformal to a fixed metric]
\label{ej:metricas-conformes-banach}
Let $(M^n,\mathbf g_0)$ be a nonempty compact Riemannian manifold
without boundary, with $n\geq1$, and let $k\in\mathbb N_0$. The set
\[
\mathcal C^k(\mathbf g_0)
:=\{a\mathbf g_0\mid a\in C^k(M,\mathbb R),\ a>0\}
\]
is a smooth embedded submanifold of the Banach space
$C^k(M,S^2T^*M)$. The map
\[
\mathcal F\colon C^k(M,\mathbb R)\longrightarrow\mathcal C^k(\mathbf g_0),
\qquad \mathcal F(f):=e^{2f}\mathbf g_0,
\]
is a smooth global parametrization with smooth inverse. If
$\mathbf g=e^{2f}\mathbf g_0$, then
\[
D\mathcal F(f)h=2h\mathbf g,
\qquad
T_{\mathbf g}\mathcal C^k(\mathbf g_0)
\simeq\{h\mathbf g\mid h\in C^k(M,\mathbb R)\}.
\]
\end{example}

\begin{proof}
Let $\mathcal A=C^k(M,\mathbb R)$ have the norm of
Lemma~\ref{lem:algebra-ck-compacta}, and let
$\mathcal B=C^k(M,S^2T^*M)$ have the analogous norm, constructed
using $\mathbf g_0$ and its Levi--Civita connection. Both spaces are
complete by that lemma. Consider the operators
\[
J\colon\mathcal A\to\mathcal B,\quad J(a)=a\mathbf g_0,
\qquad
R\colon\mathcal B\to\mathcal A,\quad
R(\mathbf T)=\frac{1}{n}\operatorname{tr}_{\mathbf g_0}\mathbf T.
\]
They are continuous and linear. Indeed, $\nabla\mathbf g_0=0$ and
$|\mathbf g_0|_{\mathbf g_0}=\sqrt n$, so
$\|J(a)\|_{\mathcal B}=\sqrt n\|a\|_{\mathcal A}$. The contraction
defining the trace is parallel and has norm $\sqrt n$ on tensors of
order two. Thus
\[
\|R(\mathbf T)\|_{\mathcal A}
\leq\frac{1}{\sqrt n}\|\mathbf T\|_{\mathcal B}.
\]
Since $R\circ J=I_{\mathcal A}$, the operator
$Q:=I_{\mathcal B}-JR$ is a continuous projection onto the closed
subspace $\mathcal B_0:=\ker R$. Each tensor has the unique
decomposition
\[
\mathbf T=J(R\mathbf T)+Q\mathbf T,
\qquad
\mathcal B=J(\mathcal A)\oplus\mathcal B_0.
\]
The image $J(\mathcal A)$ is closed: if $J(a_\nu)\to\mathbf T$,
continuity of $R$ gives $a_\nu\to R\mathbf T$, and hence
$\mathbf T=J(R\mathbf T)$.

Positivity of $a$ is an open condition on $\mathcal A$. Therefore
$\mathcal C^k(\mathbf g_0)=J(\mathcal A_+)$ is an open subset of
the complemented subspace $J(\mathcal A)$. Explicitly, if
$\mathbf g=a\mathbf g_0$, the affine map
\[
\mathbf T\longmapsto\bigl(R\mathbf T-a,\ Q\mathbf T\bigr)
\]
is an isomorphism between $\mathcal B$ and
$\mathcal A\times\mathcal B_0$, with inverse
$(b,\mathbf S)\mapsto(a+b)\mathbf g_0+\mathbf S$.
It sends $\mathbf g$ to the origin. Restricted to the open set
$R^{-1}(\mathcal A_+)$, it straightens the set of conformal metrics
onto $\mathcal A\times\{0\}$. This verifies the definition of an
embedded submanifold without using a previously constructed structure
on the space of all metrics.

By Example~\ref{ej:funciones-elementales-algebras},
$f\mapsto e^{2f}$ is a smooth diffeomorphism between $\mathcal A$
and $\mathcal A_+$. Hence
$\mathcal F=J\circ\operatorname{Exp}\circ(2I)$ is a smooth global
parametrization whose inverse is
\[
\mathbf g\longmapsto
\frac{1}{2}\log\left(\frac{1}{n}
\operatorname{tr}_{\mathbf g_0}\mathbf g\right).
\]
The chain rule gives $D\mathcal F(f)h=2he^{2f}\mathbf g_0$.
The image of this differential is $J(\mathcal A)$: multiplication
by $2e^{2f}$ is an isomorphism of $\mathcal A$ whose inverse is
multiplication by $\displaystyle\frac{1}{2}e^{-2f}$. This proves
the tangent space formula under inclusion into $\mathcal B$.
The curve $t\mapsto e^{2(f+th)}\mathbf g_0$ explicitly realizes
each velocity $2h\mathbf g$.

The same parametrization allows us to compute a geometric variation.
In coordinates, the matrix of $\mathbf g=e^{2f}\mathbf g_0$ has
determinant $e^{2nf}\det(\mathbf g_0)$, since each of its $n$ rows
is multiplied by $e^{2f}$. By the definition of Riemannian measure,
\[
d\lambda_{\mathbf g}=e^{nf}\,d\lambda_{\mathbf g_0}.
\]
Consequently, total volume expressed in the coordinate $f$ is the
smooth functional
\[
\mathcal V(f):=\int_Me^{nf}\,d\lambda_{\mathbf g_0},
\qquad
D\mathcal V(f)h=n\int_Mhe^{nf}\,d\lambda_{\mathbf g_0}.
\]
Differentiation is justified in the $C^k$ norm by smoothness of the
exponential and continuity of the integral functional. The equation
$D\mathcal V(f)h=0$ describes conformal variations preserving volume
to first order.
\end{proof}
\subsection{Curves on a sphere with fixed endpoints}

\begin{example}[A manifold of Sobolev curves]
\label{ej:curvas-esfera-extremos-banach}
Let $d\geq2$ and $p,q\in S^{d-1}\subseteq\mathbb R^d$. The set
\[
\mathcal P_{p,q}:=
\{\gamma\in H^1((0,1),\mathbb R^d)\mid
|\gamma(t)|=1\text{ for every }t\in[0,1],\
\gamma(0)=p,\ \gamma(1)=q\}
\]
is a smooth embedded submanifold of the affine space of $H^1$ curves
with these endpoints. Each Sobolev class is evaluated using its unique
continuous representative. Its tangent space at $\gamma$, under inclusion
into $H^1((0,1),\mathbb R^d)$, is
\[
Y_\gamma:=
\{h\in H^1((0,1),\mathbb R^d)\mid
h(0)=h(1)=0,\quad\gamma(t)\cdot h(t)=0\text{ for every }t\in[0,1]\}.
\]
An explicit inverse chart around $\gamma$ is
\[
\Theta_\gamma(h)(t):=
\frac{\gamma(t)+h(t)}{\sqrt{1+|h(t)|^2}},
\qquad h\in Y_\gamma.
\]
\end{example}

\begin{proof}
Let $\mathcal A_1$ be the algebra of
Lemma~\ref{lem:algebra-h1-intervalo}. On
$\mathcal H:=H^1((0,1),\mathbb R^d)$ use the equivalent norm
$\|u\|_{\mathcal H}:=\displaystyle\sum_{j=1}^d\|u_j\|_{\mathcal A_1}$.
The estimates in that lemma imply
\[
\begin{aligned}
\|u\|_{C^0([0,1],\mathbb R^d)}&\leq\|u\|_{\mathcal H},\\
\|u\cdot v\|_{\mathcal A_1}&\leq\|u\|_{\mathcal H}\|v\|_{\mathcal H},\\
\|au\|_{\mathcal H}&\leq\|a\|_{\mathcal A_1}\|u\|_{\mathcal H}.
\end{aligned}
\]
Thus the inner product and multiplication of a curve by a function
are continuous bilinear operations.

Evaluation $u\mapsto(u(0),u(1))$ is continuous, linear, and surjective.
It has the continuous right inverse
\[
(v,w)\longmapsto\bigl(t\mapsto(1-t)v+tw\bigr).
\]
Hence $\mathcal H_0:=\{h\in\mathcal H\mid h(0)=h(1)=0\}$ is
closed and complemented. The affine space
$\mathcal E_{p,q}:=\{u\in\mathcal H\mid u(0)=p,\ u(1)=q\}$
is a Banach manifold modeled on $\mathcal H_0$.
The set $\mathcal P_{p,q}$ is nonempty. If $q\neq-p$, the curve
\[
t\longmapsto\frac{(1-t)p+tq}{|(1-t)p+tq|}
\]
is smooth with the prescribed endpoints: its denominator never
vanishes, since a nontrivial convex combination of two unit vectors
can be zero only when they are opposite. If $q=-p$, choose a unit
vector $r$ perpendicular to $p$, which exists because $d\geq2$, and
concatenate the preceding curves from $p$ to $r$ and from $r$ to $q$,
reparametrized on the two halves of the interval. The resulting curve
is continuous and piecewise smooth with bounded derivative on each
half; continuity at the joining point cancels the boundary term in
integration by parts and ensures that the curve belongs to $H^1$.

To verify the submanifold property, set
$\mathcal A_{1,0}:=\{a\in\mathcal A_1\mid a(0)=a(1)=0\}$,
which is a Banach space by continuity of evaluation. The map
\[
F\colon\mathcal E_{p,q}\to\mathcal A_{1,0},
\qquad F(u):=u\cdot u-1,
\]
is well defined because $|p|=|q|=1$. It is smooth as a continuous
quadratic map. Its differential in a direction
$h\in\mathcal H_0$ is $DF(\gamma)h=2\gamma\cdot h$.
For $\gamma\in\mathcal P_{p,q}$, the operator
\[
R_\gamma\colon\mathcal A_{1,0}\to\mathcal H_0,
\qquad R_\gamma(a):=\frac{1}{2}a\gamma,
\]
is continuous and linear and satisfies $DF(\gamma)R_\gamma(a)=a$.
Consequently, zero is a split regular value of $F$.
Theorem~\ref{teo:preimagen-valor-regular-split-banach} shows that
$\mathcal P_{p,q}=F^{-1}(0)$ is an embedded submanifold and that
its tangent space is $Y_\gamma=\ker DF(\gamma)$. The projection
onto this kernel is explicitly
\[
P_\gamma h=h-(\gamma\cdot h)\gamma,
\qquad
\mathcal H_0=Y_\gamma\oplus R_\gamma(\mathcal A_{1,0}).
\]
Indeed, $\gamma\cdot P_\gamma h=0$ and
$h=P_\gamma h+R_\gamma(2\gamma\cdot h)$; the intersection of the
summands is zero because $DF(\gamma)R_\gamma=I$.

We now construct the charts with an explicit parametrization.
If $h\in Y_\gamma$, then $1+|h|^2$ belongs to $\mathcal A_1$ and
is strictly positive. Example~\ref{ej:funciones-elementales-algebras}
shows that its power $-\displaystyle\frac{1}{2}$ belongs to
$\mathcal A_1$ and depends smoothly on $h$. Hence $\Theta_\gamma$
is smooth with values in $\mathcal H$. Since $|\gamma|=1$ and
$\gamma\cdot h=0$,
\[
|\gamma+h|^2=1+|h|^2,
\qquad
\gamma\cdot\Theta_\gamma(h)=\frac{1}{\sqrt{1+|h|^2}}>0.
\]
The endpoints of $h$ are zero, so
$\Theta_\gamma(h)\in\mathcal P_{p,q}$.

Consider the relatively open set
\[
U_\gamma:=\{u\in\mathcal P_{p,q}\mid
\gamma(t)\cdot u(t)>0\text{ for every }t\in[0,1]\}.
\]
It is open because $u\mapsto\gamma\cdot u$ is continuous with values
in $C^0([0,1])$ and strict positivity on a compact set is open.
For $u\in U_\gamma$, define
\[
\varphi_\gamma(u):=\frac{u}{\gamma\cdot u}-\gamma.
\]
The denominator has positive minimum and its reciprocal belongs to
$\mathcal A_1$. Thus $\varphi_\gamma(u)\in\mathcal H_0$.
Moreover, $\gamma\cdot\varphi_\gamma(u)=0$, so it takes values
in $Y_\gamma$. The identities
\[
\varphi_\gamma(\Theta_\gamma(h))=h,
\qquad
1+|\varphi_\gamma(u)|^2=\frac{1}{(\gamma\cdot u)^2}
\]
give $\Theta_\gamma(\varphi_\gamma(u))=u$, using
$\gamma\cdot u>0$ to take the positive square root.
Both maps are continuous by the product estimates and smoothness
of the reciprocal. Therefore
$\varphi_\gamma\colon U_\gamma\to Y_\gamma$ is a global chart
on this open set. Its formulas are smooth as maps into the ambient
spaces; Lemma~\ref{lem:corestriccion-subvariedad-banach} ensures
that they are smooth also as maps into the submanifold and the
corresponding closed subspace. This chart is compatible with the
structure obtained through $F$.

On an overlap with the chart centered at $\sigma$, the coordinate
change is
\[
\varphi_\sigma\circ\Theta_\gamma(h)
=\frac{\gamma+h}{\sigma\cdot(\gamma+h)}-\sigma.
\]
Its domain is the open set where the denominator is strictly positive.
The formula and smoothness of products and reciprocals verify the
regularity of the transitions directly.

Finally, the differential of $\Theta_\gamma$ at zero is the inclusion
$h\mapsto h$. An explicit error bound can also be obtained.
If $\|h\|_{\mathcal H}<1$, then
$\||h|^2\|_{\mathcal A_1}\leq\|h\|_{\mathcal H}^2<1$.
The nonconstant coefficients in the series for
$(1+w)^{-\displaystyle\frac{1}{2}}$ have absolute value at most one because
$\binom{2n}{n}\leq4^n$. Consequently,
\[
\|\Theta_\gamma(h)-\gamma-h\|_{\mathcal H}
\leq
\frac{(\|\gamma\|_{\mathcal H}+\|h\|_{\mathcal H})
\|h\|_{\mathcal H}^2}{1-\|h\|_{\mathcal H}^2}.
\]
After division by $\|h\|_{\mathcal H}$, the right-hand side tends
to zero. The curve $s\mapsto\Theta_\gamma(sh)$ realizes each
tangent direction $h$, whereas the unit norm and endpoint constraints
force every admissible velocity to belong to $Y_\gamma$.
\end{proof}
\subsection{The Grassmannian of a Hilbert space}

\begin{example}[Subspaces as points of a manifold]
\label{ej:grassmanniana-hilbert-autocontenida}
Let $H$ be a real Hilbert space. Denote the set of its closed linear
subspaces by $\operatorname{Gr}(H)$. For
$Z\in\operatorname{Gr}(H)$, let $P_Z$ be the orthogonal projection
onto $Z$. The distance
\[
d(Z,W):=\|P_Z-P_W\|_{\mathcal L(H)}
\]
defines a Hausdorff topology on $\operatorname{Gr}(H)$.
With this topology, $\operatorname{Gr}(H)$ is a smooth Banach
manifold. Around $Y\in\operatorname{Gr}(H)$, its coordinates are
operators $A\in\mathcal L(Y,Y^\perp)$, and the corresponding
subspace is the graph
\[
\Gamma_Y(A):=\{y+Ay\mid y\in Y\}.
\]
In particular, the tangent space at $Y$ is identified through this
chart with $\mathcal L(Y,Y^\perp)$.
\end{example}

\begin{proof}
The orthogonal decomposition in
Corollary~\ref{cor:descomposicion-ortogonal-hilbert} provides $P_Z$.
If $P_Z=P_W$, their images agree and $Z=W$; the remaining properties
of $d$ follow from the operator norm. Thus $d$ is a distance.

Fix $Y$ and write $H=Y\oplus Y^\perp$. The space
$\mathcal L(Y,Y^\perp)$ is Banach by
Theorem~\ref{teo:completitud-espacio-operadores}, since closed
subspaces of a Hilbert space are complete. Let $\mathcal U_Y$ be the
set of subspaces $Z$ for which
$P_Y\restriction_Z\colon Z\to Y$ is a continuous linear
isomorphism. Each graph $\Gamma_Y(A)$ belongs to $\mathcal U_Y$
and is closed: if $y_\nu+Ay_\nu\to x$, applying $P_Y$ gives
$y_\nu\to P_Yx$, and applying $P_{Y^\perp}$ gives
$P_{Y^\perp}x=A(P_Yx)$. Conversely, if $Z\in\mathcal U_Y$, then
\[
A_Z:=P_{Y^\perp}\circ(P_Y\restriction_Z)^{-1}
\in\mathcal L(Y,Y^\perp),
\qquad Z=\Gamma_Y(A_Z).
\]
The representation is unique because the $Y$ component of $y+Ay$
is $y$.

We compute the orthogonal projection onto a graph to verify
continuity of these coordinates. Define
\[
J_A\colon Y\to H,\quad J_Ay=y+Ay,
\qquad S_A:=I_Y+A^*A.
\]
The operator $S_A$ is self-adjoint and
$\langle S_Ay,y\rangle=\|y\|^2+\|Ay\|^2$. In particular,
$\|S_Ay\|\geq\|y\|$. This inequality shows that its image is
closed: a Cauchy sequence of images has Cauchy preimages.
Moreover, a vector perpendicular to its image belongs to
$\ker S_A^*=\ker S_A=\{0\}$. The image is therefore dense and
closed, so $S_A$ is invertible with $\|S_A^{-1}\|\leq1$.
If $Y=\{0\}$, these identities are interpreted in the zero space
and the conclusion is immediate. Since $J_A^*J_A=S_A$, the operator
\[
P_{\Gamma_Y(A)}=J_AS_A^{-1}J_A^*
\]
is self-adjoint and idempotent with image $\Gamma_Y(A)$: its image
is contained in that graph, and its composition with $J_A$ is $J_A$.
It is therefore the orthogonal projection onto the graph.
By Theorem~\ref{teo:adjunto-operador-acotado-hilbert}, taking
the adjoint is continuous and linear in operator norm. Products are
continuous bilinear operations, and inversion is smooth
by Proposition~\ref{prop:inversion-operadores-suave}.
Consequently, $A\mapsto P_{\Gamma_Y(A)}$ is smooth and, in
particular, $\Gamma_Y$ is continuous for the distance $d$.

We still need to show that $\mathcal U_Y$ is open and that
$Z\mapsto A_Z$ is continuous. Take
$Z_0=\Gamma_Y(A_0)\in\mathcal U_Y$ and set $P_0:=P_{Z_0}$.
For another subspace $Z$, with $P:=P_Z$, define
\[
T:=PP_0+(I-P)(I-P_0)
=I+(P-P_0)(2P_0-I).
\]
Orthogonal decomposition gives $\|2P_0-I\|\leq1$, and
$\|T-I\|\leq\|P-P_0\|$. If $d(Z,Z_0)<1$, the Neumann series
in Theorem~\ref{teo:serie-de-neumann} shows that $T$ is invertible.
Moreover, $TP_0=PT$, so $T(Z_0)=Z$: one inclusion follows from
this identity and the other from $P_0T^{-1}=T^{-1}P$.

The operator $J_{A_0}$ is an isomorphism from $Y$ onto $Z_0$.
If
\[
d(Z,Z_0)<\delta_0,
\qquad
\delta_0:=\frac{1}{2(1+\|J_{A_0}\|)},
\]
then $\|T-I\|<1$ and
\[
\|P_YTJ_{A_0}-I_Y\|
\leq\|T-I\|\|J_{A_0}\|<\frac{1}{2}.
\]
A second use of the Neumann series proves that $P_YTJ_{A_0}$ is
invertible. Since $TJ_{A_0}$ is an isomorphism from $Y$ onto $Z$,
it follows that $P_Y\restriction_Z$ is an isomorphism; hence
$Z\in\mathcal U_Y$. In this neighborhood,
\[
A_Z=P_{Y^\perp}TJ_{A_0}(P_YTJ_{A_0})^{-1}.
\]
The right-hand side depends continuously on $P_Z$ in operator norm.
This proves openness of $\mathcal U_Y$ and continuity of the chart
$\varphi_Y\colon\mathcal U_Y\to\mathcal L(Y,Y^\perp)$,
$\varphi_Y(Z)=A_Z$. Together with continuity of $\Gamma_Y$,
it shows that these maps are inverse homeomorphisms.
Each subspace $Y$ belongs to $\mathcal U_Y$ and corresponds to the
zero operator, so these charts cover $\operatorname{Gr}(H)$.

Now let $Y,W\in\operatorname{Gr}(H)$. If
$\Gamma_Y(A)\in\mathcal U_W$, the operator
$C_A:=P_WJ_A\colon Y\to W$ is invertible. The coordinate change is
\[
\varphi_W\circ\Gamma_Y(A)
=P_{W^\perp}J_A(P_WJ_A)^{-1}.
\]
Its domain is open because $A\mapsto P_WJ_A$ is continuous and
the set of isomorphisms between two Banach spaces is open.
The same inversion proposition and the chain rule show that the
transition is smooth. Interchanging $Y$ and $W$ gives smoothness
of its inverse. By the definition of a Banach manifold, this
constructs the asserted structure without presupposing a
differentiable structure on the set of subspaces.

The chart centered at $Y$ identifies the kinematic tangent space
with $\mathcal L(Y,Y^\perp)$. To describe its realization in
operators on $H$ as well, let $B\in\mathcal L(Y,Y^\perp)$ and
consider the curve $s\mapsto\Gamma_Y(sB)$. Relative to
$H=Y\oplus Y^\perp$,
\[
P_{\Gamma_Y(sB)}=
\begin{pmatrix}
S_s^{-1}&sS_s^{-1}B^*\\
sBS_s^{-1}&s^2BS_s^{-1}B^*
\end{pmatrix},
\qquad S_s:=I_Y+s^2B^*B.
\]
Since $\|S_s^{-1}\|\leq1$,
$\|S_s^{-1}-I_Y\|\leq s^2\|B\|^2$.
The norm of a matrix of operators is at most the sum of the norms
of its four blocks. Therefore, if $|s|\|B\|\leq1$,
\[
\left\|P_{\Gamma_Y(sB)}-P_Y
-s\begin{pmatrix}0&B^*\\B&0\end{pmatrix}\right\|
\leq2s^2\|B\|^2+2|s|^3\|B\|^3
\leq4s^2\|B\|^2.
\]
For $s\neq0$, dividing by $|s|$ gives a bound tending to zero.
Thus the velocity of the projection is exactly
$\begin{pmatrix}0&B^*\\B&0\end{pmatrix}$.
This distinguishes the operator $B$ representing the tangent vector
from the projection $P_Y$ representing the point of the Grassmannian.
\end{proof}
\section{Bump functions and smooth partitions of unity}
\label{sec:particiones-unidad-banach}

Constructing global objects from local data requires additional
care on Banach manifolds. In finite dimensions, second countability
and local compactness give paracompactness, and
Euclidean bump functions yield smooth partitions of unity. In an
infinite-dimensional Banach space, there are no relatively compact
neighborhoods, and existence of bump functions with the prescribed regularity
depends on the geometry of the model space. It is therefore useful to separate
the topological and differential parts of the argument.

Let $E$ be a topological space. Recall that the support of a
continuous function $b\colon E\longrightarrow\mathbb R$ is the
closed set
\[
\supp(b):=\overline{\{x\in E\mid b(x)\neq0\}}.
\]
The closure is taken in $E$.

\begin{definition}\label{def:funcion-flan-banach}\glsadd{funcion-flan-banach}
\index{bump function!on a Banach space}
Let $E$ be a real Banach space and let
$k\in\mathbb N\cup\{\infty\}$. A \textbf{bump function of class $C^k$}
on $E$, also called a \textbf{bump of class $C^k$}, is a function
$b\in C^k(E,\mathbb R)$ that is not identically zero and has bounded support.
\end{definition}

In finite dimensions, a closed bounded support is compact, so the
preceding definition agrees with the one used in earlier chapters.
In infinite dimensions, compact support cannot be required, as the
following result shows.

\begin{proposition}\label{prop:no-hay-funciones-soporte-compacto-banach}
Let $E$ be an infinite-dimensional normed space. Every continuous function
$b\colon E\longrightarrow\mathbb R$ with compact support is identically zero.
\end{proposition}

\begin{proof}
We first use the following elementary consequence of Riesz's lemma. If
$Y\subsetneq E$ is a closed subspace and $0<a<1$, there exists $x\in E$ such that
$\|x\|_E=1$ and
\[
\operatorname{dist}(x,Y)>a.
\]
Indeed, take $z\notin Y$ and write
$d:=\operatorname{dist}(z,Y)>0$. By the definition of the infimum, there exists
$y_0\in Y$ such that $\|z-y_0\|_E<\frac{d}{a}$. If
$x:=\frac{z-y_0}{\|z-y_0\|_E}$, then, for every $y\in Y$,
\[
\|x-y\|_E
=
\frac{\|z-(y_0+\|z-y_0\|_Ey)\|_E}{\|z-y_0\|_E}
\geq
\frac{d}{\|z-y_0\|_E}
>a.
\]

Applying this observation inductively with $a=\frac{1}{2}$ and
$Y_m:=\operatorname{span}\{x_1,\dots,x_m\}$ gives a sequence
$(x_m)_{m\in\mathbb N}$ of unit vectors such that
$\|x_m-x_j\|_E>\frac{1}{2}$ whenever $m\neq j$. Here each $Y_m$ is closed because it is
finite-dimensional and proper because $E$ is infinite-dimensional. Thus
the closed unit ball of $E$ is not compact.

Now suppose $b$ is not identically zero and choose $x_0\in E$ with
$b(x_0)\neq0$. Continuity gives $r>0$ such that
\[
|b(x)|>\frac{|b(x_0)|}{2}
\qquad
\text{if }\|x-x_0\|_E<r.
\]
Consequently,
$\overline B_E(x_0,r)\subseteq\supp(b)$. If this support were compact,
so would be the closed ball $\overline B_E(x_0,r)$. Translation by
$-x_0$ and dilation by the factor $\frac{1}{r}$ would then imply compactness of the
closed unit ball, contradicting the preceding paragraph.
\end{proof}

The usual normalization of a bump function remains valid with the
notion of bounded support.

\begin{lemma}\label{lem:normalizacion-funcion-flan-banach}
If $E$ admits a bump function of class $C^k$, it admits a function
$\beta\in C^k(E,\mathbb R)$ such that
\[
0\leq\beta\leq1,
\qquad
\beta=1\text{ on a neighborhood of }0,
\qquad
\supp(\beta)\text{ is bounded}.
\]
\end{lemma}

\begin{proof}
Let $b$ be a bump function and choose $z\in E$ such that $c:=b(z)\neq0$. There exist
disjoint open intervals containing $0$ and $c$, respectively.
By Proposition~\ref{funciones flan}, applied in $\mathbb R$, we can
choose $\vartheta\in C^\infty(\mathbb R)$ such that
$0\leq\vartheta\leq1$, $\vartheta=0$ on a neighborhood of $0$, and
$\vartheta=1$ on a neighborhood of $c$. Define
\[
\beta(x):=\vartheta(b(x+z)).
\]
The function $\beta$ is of class $C^k$ by the chain rule, equals one on
a neighborhood of $0$, and takes values in $[0,1]$. Moreover,
\[
\{x\in E\mid\beta(x)\neq0\}
\subseteq
\supp(b)-z,
\]
because $\vartheta(0)=0$. Since the set on the right is closed and
bounded, $\supp(\beta)$ is also bounded.
\end{proof}

Differentiable norms are a direct source of bump functions.

\begin{proposition}\label{prop:norma-suave-produce-flan-banach}
Let $\nu$ be a norm on a Banach space $E$ and suppose there exist
constants $c_\nu,C_\nu>0$ such that
\[
 c_\nu\|x\|_E\leq\nu(x)\leq C_\nu\|x\|_E,
 \qquad x\in E.
\]
If $\nu\colon E\setminus\{0\}\longrightarrow\mathbb R$ is of class $C^k$,
then
$E$ admits a bump function of class $C^k$.
\end{proposition}

\begin{proof}
Choose $\chi\in C^\infty(\mathbb R)$ such that $0\leq\chi\leq1$,
$\chi(t)=1$ for $t\leq1$, and $\chi(t)=0$ for $t\geq2$. The function
\[
b(x):=\chi(\nu(x))
\]
is of class $C^k$ away from the origin and constant on the neighborhood
$\{x\in E\mid\nu(x)<1\}$ of $0$; hence it is of class $C^k$ on all of $E$.
Moreover, $b(0)=1$ and
$\supp(b)\subseteq\{x\in E\mid\nu(x)\leq2\}$; the inequality
$c_\nu\|x\|_E\leq\nu(x)$ shows that this support is contained in
$\left\{x\in E\mathrel{\big|}\|x\|_E\leq\frac{2}{c_\nu}\right\}$.
\end{proof}

Existence of a bump function is necessary for constructing smooth
partitions of unity, but the converse implication will not be used without
additional assumptions. The precise condition needed in the geometric
argument is formulated next.

\begin{definition}\label{def:modelo-admite-particiones-Ck}
\index{partition of unity!on a Banach space}
Let $E$ be a Banach space and let
$k\in\mathbb N\cup\{\infty\}$. We say that $E$ \textbf{admits partitions
of unity of class $C^k$} if every open cover
$\mathscr U=(U_\alpha)_{\alpha\in A}$ of $E$ admits a partition of
unity $(\psi_\alpha)_{\alpha\in A}$ subordinate to $\mathscr U$, in the
sense of
Definition~\ref{def:nociones-fundamentales-espacio-topologico-cubierta-abierta-particion-unidad},
such that each $\psi_\alpha$ is of class $C^k$.
\end{definition}

This property includes, in particular, the ability to separate a closed set
from an open set by a smooth function.

\begin{lemma}[Cutoff functions in the model space]
\label{lem:cortes-modelo-particiones-banach}
Suppose $E$ admits partitions of unity of class $C^k$. If
$A\subseteq E$ is closed, $B\subseteq E$ is open, and $A\subseteq B$, there exists
$\theta\in C^k(E,\mathbb R)$ such that
\[
0\leq\theta\leq1,
\qquad
\theta=1\text{ on }A,
\qquad
\supp(\theta)\subseteq B.
\]
\end{lemma}

\begin{proof}
Consider the open cover $(B,E\setminus A)$ of $E$ and a subordinate partition
of unity $(\theta,\zeta)$ of class $C^k$. Subordination
gives $\supp(\theta)\subseteq B$ and
$\supp(\zeta)\subseteq E\setminus A$. In particular, $\zeta=0$ on $A$;
since $\theta+\zeta=1$, it follows that $\theta=1$ on $A$.
\end{proof}

Taking $A=\overline B_E(0,1)$ and $B=B_E(0,2)$ in the preceding lemma gives
a bump function of class $C^k$. Thus admitting $C^k$ partitions is,
a priori, a stronger condition than having a single
bump function.

\subsection{The global construction on a Banach manifold}

The topological part of the construction is collected in a shrinking lemma. The
paracompactness in its statement is paracompactness of the underlying
topological space; it is not part of
Definition~\ref{def:variedad-banach} and must therefore be imposed
explicitly.

\begin{lemma}[Locally finite shrinking]
\label{lem:encogimiento-localmente-finito-banach}
Let $X$ be a paracompact Hausdorff space. Then:
\begin{enumerate}[label=(\alph*)]
\item if $F\subseteq U\subseteq X$, where $F$ is closed and $U$ is open,
there exists an open set $V$ such that
\[
F\subseteq V\subseteq\overline V\subseteq U;
\]
\item if $(W_i)_{i\in I}$ is a locally finite open cover of $X$,
there exists an open cover $(V_i)_{i\in I}$ such that
\[
\overline{V_i}\subseteq W_i
\qquad\text{for every }i\in I.
\]
\end{enumerate}
\end{lemma}

\begin{proof}
First prove that $X$ is normal. Fix a point $x\in X$ and a closed set
$F$ not containing it. For each $y\in F$, the Hausdorff property
gives disjoint open sets $A_y$ and $B_y$ with $x\in A_y$ and $y\in B_y$.
The cover
\[
\{X\setminus F\}\cup\{B_y\mid y\in F\}
\]
admits a locally finite open refinement $(R_j)_{j\in J}$. Let
$J_F:=\{j\in J\mid R_j\cap F\neq\varnothing\}$. If $j\in J_F$, the set
$R_j$ cannot be contained in $X\setminus F$; since this is a refinement,
it is contained in some $B_y$. Consequently,
$x\notin\overline{R_j}$. The family $(\overline{R_j})_{j\in J}$ is also
locally finite, and hence
\[
\overline{\bigcup_{j\in J_F}R_j}
\subseteq
\bigcup_{j\in J_F}\overline{R_j}.
\]
The open set $\displaystyle W:=\bigcup_{j\in J_F}R_j$ contains $F$, while
$x\notin\overline W$. Thus $X\setminus\overline W$ is a neighborhood of
$x$ whose closure is disjoint from $F$. This proves that $X$ is regular.

Now let $F_0$ and $F_1$ be disjoint closed sets. For each $x\in F_0$,
regularity gives an open set $U_x$ such that
$x\in U_x$ and $\overline{U_x}\cap F_1=\varnothing$. Take a locally finite refinement of
the cover
\[
\{X\setminus F_0\}\cup\{U_x\mid x\in F_0\}
\]
and let $G$ be the union of those elements of the refinement that meet
$F_0$. These elements are contained in some $U_x$. Local finiteness
implies that the closure of their union is contained in the union of their
closures; consequently,
\[
F_0\subseteq G,
\qquad
\overline G\cap F_1=\varnothing.
\]
The open sets $G$ and $X\setminus\overline G$ separate $F_0$ and $F_1$. Thus
$X$ is normal. Applying this separation to the closed sets
$F$ and $X\setminus U$ gives (a).

To prove (b), well-order $I$ and construct the open sets $V_i$ by
transfinite recursion. Suppose the $V_j$ with $j<i$ have been chosen so that
the family
\[
\{V_j\mid j<i\}\cup\{W_j\mid j\geq i\}
\]
covers $X$. Define
\[
K_i
:=
X\setminus
\left(
\bigcup_{j<i}V_j\cup\bigcup_{j>i}W_j
\right).
\]
The set $K_i$ is closed, and the inductive property shows that
$K_i\subseteq W_i$. By (a), there exists an open set $V_i$ such that
$K_i\subseteq V_i$ and $\overline{V_i}\subseteq W_i$. This choice preserves
the covering property at the next successor stage.

It also preserves this property at every limit ordinal. Indeed, if a point $x$
belonged neither to any of the already chosen $V_j$ nor to any of the remaining $W_j$,
local finiteness would allow us to take the largest index, among
those strictly preceding the limit ordinal, of the sets $W_j$
containing $x$. The covering property immediately after that
index would force $x$ to belong to one of the already chosen $V_j$, a
contradiction. At the end of the recursion, $(V_i)_{i\in I}$ covers
$X$ and satisfies $\overline{V_i}\subseteq W_i$. Since $V_i\subseteq W_i$, the
family $(V_i)_{i\in I}$ is locally finite.
\end{proof}

In metric spaces, paracompactness follows directly
from the distance, without assumptions of separability or local compactness.

\begin{theorem}[Paracompactness of metric spaces]
\label{teo:paracompacidad-espacios-metricos}
\index{paracompactness!of metric spaces}
Every open cover of a metric space admits a locally finite open
refinement. In particular, every metrizable space and every
subspace of it are paracompact.
\end{theorem}

\begin{proof}
Let $(X,d)$ be a metric space and let $(U_\alpha)_{\alpha\in A}$ be an
open cover. If $X=\varnothing$, the assertion is immediate.
Suppose $X\neq\varnothing$ and choose a well-order on $A$.
For each $\alpha$, define
\[
 f_\alpha(x):=\min\{1,d(x,X\setminus U_\alpha)\},
 \qquad x\in X,
\]
and take $f_\alpha\equiv1$ when $U_\alpha=X$. Distance to a
nonempty set is $1$--Lipschitz: from
$d(x,z)\leq d(x,y)+d(y,z)$, taking infima over $z$ in the set
gives one inequality, and the other follows by interchanging
$x$ and $y$. Since truncation at $1$ preserves this property, each $f_\alpha$
is $1$--Lipschitz. Moreover, $f_\alpha(x)>0$ if and only if $x\in U_\alpha$,
since $U_\alpha$ is open.

For $n\in\mathbb N$, set $r_n:=2^{-n}$ and consider the closed sets
\[
 F_{\alpha,n}:=
 \{x\in X\mid f_\alpha(x)\geq r_n,
 \ f_\beta(x)\leq r_n/2\text{ for every }\beta<\alpha\}.
\]
They are closed because they are defined as intersections of closed
sets. If $\alpha<\beta$, $x\in F_{\alpha,n}$, and
$y\in F_{\beta,n}$, then
\[
 d(x,y)\geq |f_\alpha(x)-f_\alpha(y)|\geq r_n/2.
\]
Moreover, these sets cover $X$ as $\alpha$ and $n$ vary.
Indeed, given $x\in X$, let $\alpha$ be the first index for which
$x\in U_\alpha$. Then $f_\alpha(x)>0$ and $f_\beta(x)=0$ for every
$\beta<\alpha$. Choosing $n$ with $r_n\leq f_\alpha(x)$ suffices to obtain
$x\in F_{\alpha,n}$.

If $F_{\alpha,n}\neq\varnothing$, define the open set
\[
 V_{\alpha,n}:=\{x\in X\mid d(x,F_{\alpha,n})<r_n/8\};
\];
if $F_{\alpha,n}=\varnothing$, take $V_{\alpha,n}=\varnothing$.
For $x\in V_{\alpha,n}$, there exists $u\in F_{\alpha,n}$ with
$d(x,u)<r_n/8$, and hence
$f_\alpha(x)\geq f_\alpha(u)-d(x,u)>7r_n/8>0$. Thus
$V_{\alpha,n}\subseteq U_\alpha$.
If $\alpha\neq\beta$, $x\in V_{\alpha,n}$, and
$y\in V_{\beta,n}$, choose points in the respective closed sets at
distance less than $r_n/8$. Separation of the closed sets and the
triangle inequality give $d(x,y)>r_n/4$. Consequently, each ball
$B_d(z,r_n/8)$ meets at most one of the sets
$V_{\alpha,n}$. For each fixed $n$, the family
$(V_{\alpha,n})_{\alpha\in A}$ is therefore locally finite.
Together, these families cover $X$ as $n$ varies.

To obtain local finiteness also with respect to $n$, set
\[
 O_n:=\bigcup_{1\leq k\leq n}\ \bigcup_{\alpha\in A}V_{\alpha,k},
 \qquad
 h_n(x):=\min\{1,d(x,X\setminus O_n)\},
\],
with $h_n\equiv1$ if $O_n=X$, and define
\[
 P_0:=\varnothing,
 \qquad
 P_n:=\{x\in X\mid h_n(x)\geq2^{-n}\}
 \quad(n\in\mathbb N).
\]
The sets $O_n$ are open and increasing. The functions $h_n$ are
continuous and satisfy $h_n\leq h_{n+1}$; hence the $P_n$ are
closed, $P_n\subseteq O_n$, and $P_n\subseteq P_{n+1}$. Their interiors
cover $X$. Indeed, if $x\in O_m$, then $h_m(x)>0$ and, for sufficiently large $n$,
$n\geq m$ and
$h_n(x)\geq h_m(x)>2^{-n}$. By continuity, this strict inequality
persists on a neighborhood of $x$ contained in $P_n$.

Finally, define
\[
 W_{\alpha,n}:=V_{\alpha,n}\setminus P_{n-1}.
\]
Each $W_{\alpha,n}$ is open and contained in $U_\alpha$.
To check that they cover $X$, given $x\in X$, take the least $n$ such
that $x\in V_{\alpha,n}$ for some $\alpha$. If $n=1$, we have
$x\notin P_0$. If $n>1$, minimality gives $x\notin O_{n-1}$ and, since
$P_{n-1}\subseteq O_{n-1}$, also $x\notin P_{n-1}$. In either case,
$x\in W_{\alpha,n}$.

Now fix $x\in X$ and choose $N$ with $x\in\operatorname{Int}(P_N)$.
This neighborhood meets no $W_{\alpha,n}$ with $n>N$, because
$P_N\subseteq P_{n-1}$. For each of the finitely many levels
$1\leq n\leq N$, local finiteness of
$(V_{\alpha,n})_{\alpha\in A}$ gives a neighborhood of $x$
meeting only finitely many of its elements. The intersection of
these neighborhoods with $\operatorname{Int}(P_N)$ meets only finitely
many of the $W_{\alpha,n}$. Thus the nonempty members of this
family form the required locally finite open refinement.
The assertion for metrizable spaces follows by choosing a compatible
metric, and the assertion for their subspaces follows by applying what we have proved
to the restricted metric.
\end{proof}

\begin{theorem}[Smooth partitions on Banach manifolds]
\label{teo:particiones-unidad-variedad-banach}\glsadd{particion-unidad-banach}
Let $M$ be a Banach manifold of class $C^k$ modeled on a Banach
space $E$, where $k\in\mathbb N\cup\{\infty\}$. Suppose the topological
space $M$ is paracompact and $E$ admits partitions of unity of
class $C^k$. Then every open cover
$\mathscr U=(U_\alpha)_{\alpha\in A}$ of $M$ admits a partition of
unity of class $C^k$ subordinate to $\mathscr U$.
\end{theorem}

\begin{proof}
Lemma~\ref{lem:encogimiento-localmente-finito-banach} shows, in
particular, that $M$ is regular. For each $x\in M$, choose
$\alpha(x)\in A$ such that $x\in U_{\alpha(x)}$ and a chart
$(C_x,\varphi_x)$ around $x$. Applying regularity to the open set
$C_x\cap U_{\alpha(x)}$ gives an open set $Q_x$ such that
\[
x\in Q_x,
\qquad
\overline{Q_x}\subseteq C_x\cap U_{\alpha(x)}.
\]
Since $\varphi_x(Q_x)$ is an open neighborhood of $\varphi_x(x)$, there exists
$r_x>0$ such that
\[
\overline B_E(\varphi_x(x),r_x)
\subseteq
\varphi_x(Q_x).
\]
Set
\[
P_x:=\varphi_x^{-1}\bigl(B_E(\varphi_x(x),r_x)\bigr).
\]
Since $P_x\subseteq Q_x$ and $\overline{Q_x}\subseteq C_x$, the closure
of $P_x$ in $M$ remains inside the chart domain. Thus
\begin{equation}
\overline{P_x}
=
\varphi_x^{-1}\bigl(\overline B_E(\varphi_x(x),r_x)\bigr)
\subseteq
Q_x
\subseteq
U_{\alpha(x)}.
\label{eq:nucleo-carta-particion-banach}
\end{equation}

Paracompactness gives a locally finite open refinement
$(W_i)_{i\in I}$ of the cover $(P_x)_{x\in M}$. For each $i$, fix a
point $x(i)$ such that $W_i\subseteq P_{x(i)}$ and abbreviate
\[
C_i:=C_{x(i)},
\qquad
\varphi_i:=\varphi_{x(i)},
\qquad
\alpha(i):=\alpha(x(i)).
\]
By Lemma~\ref{lem:encogimiento-localmente-finito-banach}, there exists an
open cover $(V_i)_{i\in I}$ with
\begin{equation}
\overline{V_i}\subseteq W_i
\qquad\text{for every }i\in I.
\label{eq:encogimiento-particion-banach}
\end{equation}

Within the model space, consider
\[
A_i:=\varphi_i(\overline{V_i}),
\qquad
B_i:=\varphi_i(W_i).
\]
The set $B_i$ is open and $A_i\subseteq B_i$. Moreover, $A_i$ is closed
in $E$. To verify this last point, first observe that
$\overline{V_i}$ is closed in $M$ and contained in $C_i$, so
$A_i$ is closed in the open set $\varphi_i(C_i)$. On the other hand,
$A_i$ is contained in the closed ball
$\overline B_E(\varphi_i(x(i)),r_{x(i)})$, which is a closed subset of $E$
contained in $\varphi_i(C_i)$ by \eqref{eq:nucleo-carta-particion-banach}.
If $A_i=F_i\cap\varphi_i(C_i)$ with $F_i$ closed in $E$, then
\[
A_i
=
F_i\cap
\overline B_E(\varphi_i(x(i)),r_{x(i)}),
\],
proving the assertion.

Lemma~\ref{lem:cortes-modelo-particiones-banach} gives a function
$\theta_i\in C^k(E,\mathbb R)$ such that
\begin{equation}
0\leq\theta_i\leq1,
\qquad
\theta_i=1\text{ on }A_i,
\qquad
\supp(\theta_i)\subseteq B_i.
\label{eq:corte-modelo-particion-banach}
\end{equation}
Define $\eta_i\colon M\longrightarrow\mathbb R$ by
\[
\eta_i(p)
:=
\begin{cases}
\theta_i(\varphi_i(p)),&p\in C_i,\\
0,&p\notin C_i.
\end{cases}
\]
This extension by zero is of class $C^k$. Indeed, the third property in
\eqref{eq:corte-modelo-particion-banach} implies that the local function
vanishes outside a closed set contained in $W_i$, and
$\overline{W_i}\subseteq\overline{P_{x(i)}}\subseteq C_i$ by
\eqref{eq:nucleo-carta-particion-banach}; hence it is identically zero on
a neighborhood of every point of $M\setminus C_i$. Moreover,
\begin{equation}
\supp(\eta_i)\subseteq W_i,
\qquad
\eta_i=1\text{ on }\overline{V_i}.
\label{eq:corte-variedad-particion-banach}
\end{equation}

The family $(\supp(\eta_i))_{i\in I}$ is locally finite because it is
subordinate to $(W_i)_{i\in I}$. Consequently, the sum
\[
h:=\sum_{i\in I}\eta_i
\]
locally consists of finitely many functions of class $C^k$ and
therefore belongs to $C^k(M)$. Given $p\in M$, there exists $i$ with $p\in V_i$;
then \eqref{eq:corte-variedad-particion-banach} gives $\eta_i(p)=1$. Thus
$h(p)\geq1$ for every $p\in M$. The functions
\[
\rho_i:=\frac{\eta_i}{h}
\]
are of class $C^k$, are nonnegative, form a locally finite family,
satisfy $\supp(\rho_i)\subseteq W_i$, and obey
$\displaystyle\sum_{i\in I}\rho_i=1$.

Finally, for each $\alpha\in A$, define
\[
\psi_\alpha
:=
\sum_{\{i\in I\mid\alpha(i)=\alpha\}}\rho_i.
\]
The sum is locally finite, so $\psi_\alpha$ is of class $C^k$.
Moreover,
\[
\sum_{\alpha\in A}\psi_\alpha
=
\sum_{i\in I}\rho_i
=1.
\]
To check subordination, note that any subfamily of a
locally finite family of closed sets has closed union. Hence
\[
\supp(\psi_\alpha)
\subseteq
\bigcup_{i\in I,\ \alpha(i)=\alpha}\supp(\rho_i)
\subseteq
\bigcup_{i\in I,\ \alpha(i)=\alpha}W_i
\subseteq
U_\alpha.
\]
The family of these supports is also locally finite, since it is obtained
by grouping the locally finite family $(\supp(\rho_i))_{i\in I}$ according to the
elements of the original cover. Thus
$(\psi_\alpha)_{\alpha\in A}$ is the required partition.
\end{proof}

\begin{definition}\label{def:variedad-banach-Ck-paracompacta}
\index{Banach manifold!Ck-paracompact@$C^k$--paracompact}
A Banach manifold of class $C^k$ is called
\textbf{$C^k$--paracompact} if its topological space is paracompact and every
open cover admits a subordinate partition of unity of class $C^k$.
\end{definition}

Theorem~\ref{teo:particiones-unidad-variedad-banach} states that a
paracompact manifold modeled on a space admitting partitions of
class $C^k$ is $C^k$--paracompact. This formulation distinguishes the two
assumptions: paracompactness produces locally finite families, while
the property of the model space provides differentiable cutoff functions.

\subsection{Criteria and examples of model spaces}

General criteria for passing from a single bump function to
partitions of unity belong to the theory of smoothness and renorming
of Banach spaces. We introduce the functional assumption appearing in
one of the most useful criteria.

\begin{definition}\label{def:espacio-asplund}
\index{Asplund space}
A Banach space $E$ is an \textbf{Asplund space} if every closed separable
subspace $Y\subseteq E$ has separable dual.
\end{definition}

\begin{theorem}[Fry's criterion]
\label{teo:criterio-fry-particiones-banach}
Let $E$ be a Banach space and let
$k\in\mathbb N\cup\{\infty\}$. If $E$ admits a bump function of class
$C^k$ and $E'$ is an Asplund space, then $E$ admits partitions of
unity of class $C^k$.
\end{theorem}

This result is Fry's Theorem~3.1~\cite[p.~146]{Fry1998}. Its proof uses projectional resolutions in dual spaces and constructs a topological embedding into a space $c_0(\Gamma)$ whose scalar components are of class $C^k$. Existence of this embedding provides the partitions of unity. The criterion thus combines the regularity of the bump function with the Asplund assumption on the dual. Theorem~\ref{teo:particiones-unidad-variedad-banach} explains how to transfer its conclusion to a paracompact manifold. Johanis~\cite{Johanis2017} gives a characterization that brings together and extends several existence criteria.

In that work, sufficiency of a single bump function of class $C^k$ on an arbitrary Banach space is formulated as an open problem. In addition to the bump function, his characterization uses a family of smooth functions that allows each point to be approximated by subspaces determined by its nonzero components. This additional information makes partitions possible in the cases covered by the criterion.

Recall why the Asplund assumption in the preceding criterion holds in
the reflexive examples. If $E$ is reflexive, so is $E'$ by
Proposition~\ref{prop: reflexivo iff dual reflexivo}. Let $Y\subseteq E'$ be a
closed separable subspace. Since closed subspaces of a reflexive
space are reflexive, $Y$ is reflexive: indeed, the unit ball of
$Y$ is the intersection of the unit ball of $E'$ with the closed
convex subspace $Y$, which is weakly closed by
Theorem~\ref{teo: cerradura convexa fuerte y debil}; thus this ball is
weakly compact, and Theorem~\ref{teo: kakutani} gives reflexivity of
$Y$. Consequently, $Y'$ is also reflexive, and its closed unit ball
$B_{Y'}$ is weakly compact.

On the other hand, separability of $Y$ and
Theorem~\ref{teo: bola dual debil estrella metrizable iff separable} show
that $B_{Y'}$ is metrizable in the weak-$*$ topology. Since $Y$ is
reflexive, the canonical embedding $J_Y\colon Y\longrightarrow Y''$, defined
by $J_Y(y)(y'):=y'(y)$, is surjective. Consequently, the weak topology
$\sigma(Y',Y'')$ and the weak-$*$ topology
$\sigma(Y',J_Y(Y))$ agree on $Y'$; through $J_Y$, the latter is
precisely $\sigma(Y',Y)$. Thus $B_{Y'}$ is a compact metric space. Choose a
countable weakly dense subset $D\subseteq B_{Y'}$. Its
real convex hull $\operatorname{co}(D)$ has weak closure
$B_{Y'}$, since it contains $D$ and is contained in the convex ball
$B_{Y'}$. The rational convex hull $\operatorname{co}_{\mathbb Q}(D)$
is countable and norm dense in $\operatorname{co}(D)$: simply approximate
the coefficients of each finite convex combination by nonnegative rational
numbers summing to one. Finally, the norm closure of
$\operatorname{co}_{\mathbb Q}(D)$ is convex and, by
Theorem~\ref{teo: cerradura convexa fuerte y debil}, also weakly
closed. Since it contains $\operatorname{co}(D)$, it contains its weak closure,
which is $B_{Y'}$; the reverse inclusion is immediate. Consequently,
$B_{Y'}$ is norm separable, and so is $Y'$. This proves that
every reflexive space is Asplund and, in particular, the dual of a
reflexive space is Asplund.

\begin{example}[Hilbert spaces]
\label{ej:particiones-unidad-hilbert}
Let $H$ be a real Hilbert space. The quadratic function
\[
Q\colon H\longrightarrow\mathbb R,
\qquad
Q(x):=\langle x,x\rangle_H,
\]
is smooth, and its only nonzero derivatives are
\[
DQ(x)h=2\langle x,h\rangle_H,
\qquad
D^2Q(x)(h_1,h_2)=2\langle h_1,h_2\rangle_H.
\]
If $\chi\in C^\infty(\mathbb R)$ satisfies $0\leq\chi\leq1$,
$\chi=1$ on $(-\infty,1]$, and $\chi=0$ on $[4,\infty)$, then
$\chi\circ Q$ is a smooth bump function.
Since $H$ is reflexive, $H'$ is Asplund.
Theorem~\ref{teo:criterio-fry-particiones-banach} implies that $H$ admits
smooth partitions of unity. Consequently, every smooth
paracompact manifold modeled on a Hilbert space admits smooth partitions
of unity, without any separability assumption on the model
space.
\end{example}

The degree of differentiability available in $L^p(X)$ spaces can be seen
directly in the functional defining their norm.

\begin{lemma}\label{lem:potencia-norma-Lp-suavidad}
Let $(X,\mathcal A,\mu)$ be a measure space and $1<p<\infty$, over real
scalars. Define
\[
P_p\colon L^p(X)\longrightarrow\mathbb R,
\qquad
P_p(u):=\int_X|u|^p\,d\mu.
\]
If $p$ is an even integer, $P_p$ is smooth. If $p$ is not an even integer, $P_p$
is of class $C^m$ for every $m\in\mathbb N$ with $m<p$. For $1\leq j<p$, its derivatives
are
\begin{equation}
D^jP_p(u)(h_1,\dots,h_j)
=
\int_X a_{p,j}(u)h_1\cdots h_j\,d\mu,
\label{eq:derivadas-potencia-norma-Lp}
\end{equation},
where
\[
a_{p,j}(t)
:=
p(p-1)\cdots(p-j+1)|t|^{p-j}(\operatorname{sgn}t)^j,
\qquad
a_{p,j}(0):=0.
\]
\end{lemma}

\begin{proof}
First suppose $1\leq j<p$. Hölder's inequality in
Proposition~\ref{desigualdad de holder}, with exponents
$\frac{p}{p-j},p,\dots,p$, gives
\begin{equation}
\left|
\int_Xa_{p,j}(u)h_1\cdots h_j\,d\mu
\right|
\leq
C_{p,j}\|u\|_{L^p(X)}^{p-j}
\prod_{\ell=1}^j\|h_\ell\|_{L^p(X)}.
\label{eq:cota-derivadas-potencia-norma-Lp}
\end{equation}
Here and in the following estimates, $C_{p,j}>0$ depends only on $p$ and
$j$; it does not depend on $u$, the increments, or the measure space.
Thus the right-hand side of
\eqref{eq:derivadas-potencia-norma-Lp} defines a continuous $j$--linear
form.

The Nemytskii operator
\begin{equation}
L^p(X)\longrightarrow L^{\frac{p}{p-j}}(X),
\qquad
u\longmapsto a_{p,j}(u),
\label{eq:nemytskii-potencia-Lp}
\end{equation}
is continuous. Indeed, if $r:=p-j$, for $0<r\leq1$ we have the pointwise
inequality
\[
\bigl||s|^r(\operatorname{sgn}s)^j
-|t|^r(\operatorname{sgn}t)^j\bigr|
\leq C_r|s-t|^r,
\],
while for $r>1$,
\[
\bigl||s|^r(\operatorname{sgn}s)^j
-|t|^r(\operatorname{sgn}t)^j\bigr|
\leq
C_r(|s|^{r-1}+|t|^{r-1})|s-t|.
\]
In both inequalities, $C_r>0$ can be chosen to depend only on
$r$; in particular, it is independent of $s$ and $t$.
Raising to the power $\frac{p}{r}$, integrating, and applying Hölder in the second case
gives, respectively,
\[
\|a_{p,j}(u)-a_{p,j}(v)\|_{L^{\frac{p}{r}}(X)}
\leq C_{p,j}\|u-v\|_{L^p(X)}^{r}
\]
and
\[
\|a_{p,j}(u)-a_{p,j}(v)\|_{L^{\frac{p}{r}}(X)}
\leq
C_{p,j}
(\|u\|_{L^p(X)}^{r-1}+\|v\|_{L^p(X)}^{r-1})
\|u-v\|_{L^p(X)}.
\]
This proves continuity of \eqref{eq:nemytskii-potencia-Lp} and, together with
\eqref{eq:cota-derivadas-potencia-norma-Lp}, norm continuity of the
function assigning to $u$ the multilinear form in
\eqref{eq:derivadas-potencia-norma-Lp}.

It remains to verify that these forms are indeed the Fréchet derivatives.
Set $a_{p,0}(t):=|t|^p$. If $j+1<p$, the fundamental theorem of calculus
on $\mathbb R$ gives, at almost every point of $X$,
\[
a_{p,j}(u+h)-a_{p,j}(u)-a_{p,j+1}(u)h
=
\int_0^1
\bigl(a_{p,j+1}(u+th)-a_{p,j+1}(u)\bigr)h\,dt.
\]
For $j=0$, interpret the empty product as $1$. For each $t\in[0,1]$, Hölder's inequality with exponents $\frac{p}{p-j-1},p,\dots,p$ controls the $L^1(X)$ norm of the integrand after multiplication by $h_1\cdots h_j$. Moreover, continuity of the Nemytskii operator makes this integrand depend continuously on $t$ in $L^1(X)$. We may therefore integrate in $t$ in this Banach space and apply the continuous linear functional $\displaystyle v\longmapsto\int_Xv\,d\mu$; this operation justifies interchanging the integrals. The preceding scalar identity agrees with this identity in $L^1(X)$: the integrand vanishes outside the union of the sets $\{|u|\geq1/n\}$ and $\{|h|\geq1/n\}$, with $n\in\mathbb N$, which have finite measure, so Fubini's theorem applies on this support of $\sigma$--finite measure. We thus obtain
\begin{align*}
&\left|
\int_X
\bigl(a_{p,j}(u+h)-a_{p,j}(u)-a_{p,j+1}(u)h\bigr)
h_1\cdots h_j\,d\mu
\right|\\
&\qquad\leq
\sup_{t\in[0,1]}
\|a_{p,j+1}(u+th)-a_{p,j+1}(u)\|_{L^{\frac{p}{p-j-1}}(X)}
\|h\|_{L^p(X)}
\prod_{\ell=1}^j\|h_\ell\|_{L^p(X)}.
\end{align*}
Continuity of the Nemytskii operator corresponding to $j+1$
shows that the supremum tends to zero as $h\to0$ in $L^p(X)$: in the two
preceding estimates, simply replace $v$ by $u+th$, and the bounds are
uniform for $0\leq t\leq1$. Thus, for every $\varepsilon>0$, there exists
$\delta>0$, depending on $u,p,j$ and $\varepsilon$, such that, if
$\|h\|_{L^p(X)}<\delta$, the norm of the remainder as a $j$--linear form is at
most $\varepsilon\|h\|_{L^p(X)}$. For $j=0$, this
proves the first derivative formula. Suppose that, for an integer $j\geq1$ with $j+1<p$, the form in \eqref{eq:derivadas-potencia-norma-Lp} already agrees with $D^jP_p$. The remainder estimate, taken over all $h_\ell\in L^p(X)$ with $\|h_\ell\|_{L^p(X)}\leq1$ and $\ell\in\{1,\ldots,j\}$, proves differentiability of $D^jP_p$ in the norm of the space of $j$--linear forms. Through the association in Proposition~\ref{prop:completitud-multilineales}, its derivative is exactly the form indicated for order $j+1$. Continuity of this form was already established in \eqref{eq:nemytskii-potencia-Lp}. This proves the induction step and hence that $P_p$ is of class $C^m$ for every $m\in\mathbb N$ with $m<p$.

If $p$ is an even integer, $|u|^p=u^p$ and
\[
(u_1,\dots,u_p)\longmapsto\int_Xu_1\cdots u_p\,d\mu
\]
is a continuous $p$--linear form by Hölder's inequality. Thus $P_p$ is a continuous homogeneous
polynomial of degree $p$. Proposition~\ref{prop:diferenciabilidad-operadores-multilineales}, applied to this form, and Corollary~\ref{cor:composicion-Ck-banach}, applied to the diagonal map $u\longmapsto(u,\ldots,u)$, show that $P_p$ is smooth.
\end{proof}

\begin{corollary}[The models $\ell^p(\mathbb N)$ and $L^p(X)$]
\label{cor:particiones-unidad-lp-banach}
Let $1<p<\infty$.
\begin{enumerate}[label=(\alph*)]
\item If $p$ is an even integer, the real spaces $\ell^p(\mathbb N)$ and
$L^p(X)$, for every $\sigma$--finite measure space
$(X,\mathcal A,\mu)$, admit smooth partitions of unity.
\item If $p$ is not an even integer, these spaces admit partitions of
unity of class $C^m$ for each $m\in\mathbb N$ with $m<p$.
\end{enumerate}
Consequently, Theorem~\ref{teo:particiones-unidad-variedad-banach}
applies, with the same degree of differentiability, to every paracompact Banach manifold
modeled on one of these spaces.
\end{corollary}

\begin{proof}
Let $\chi\in C^\infty(\mathbb R)$ be such that $0\leq\chi\leq1$, equal to one on
$(-\infty,1]$ and zero on $[2,\infty)$. By
Lemma~\ref{lem:potencia-norma-Lp-suavidad}, the function
$u\mapsto\chi(P_p(u))$ is a smooth bump function if $p$ is even and of class
$C^m$ if $m<p$ in the other cases.
Theorem~\ref{teo:dualidad-reflexividad-Lp} shows that the
$L^p(X)$ spaces under consideration are reflexive; the case of $\ell^p(\mathbb N)$ follows by taking
counting measure on $\mathbb N$. As we saw before
Example~\ref{ej:particiones-unidad-hilbert}, the dual of a reflexive space
is Asplund. The conclusion follows from
Theorem~\ref{teo:criterio-fry-particiones-banach}.
\end{proof}

The differentiability range in the corollary is not a defect of the chosen
bump function. The following result shows that, for $\ell^p(\mathbb N)$, the
restriction is intrinsic.

\begin{theorem}[Bonic--Frampton]
\label{teo:obstruccion-bonic-frampton-lp}
Let $1<p<\infty$ and suppose $p$ is not an even integer. If
$m\in\mathbb N$ and $m\geq p$, there exists no nonzero bump function on $\ell^p(\mathbb N)$
of class $C^m$. In particular, $\ell^p(\mathbb N)$ admits neither partitions of
unity of class $C^m$ for these values of $m$ nor smooth partitions of
unity.
\end{theorem}

The statement is Theorem~1 and its Corollary~1 in Bonic and Frampton
\cite[Theorem~1 and Corollary~1, pp.~393--394]{BonicFrampton1965}. Their argument
uses the following boundary value principle: if $p$ is not an even
integer, $D\subseteq\ell^p(\mathbb N)$ is open, connected, and bounded, and a real function is
continuous on $\overline D$ and $m$ times differentiable on $D$, with $m\geq p$,
then its image on $\overline D$ is contained in the closure of its image
on $\partial D$. Applied to a bounded connected component of the set where
a bump function is nonzero, this principle forces the function to be identically zero. Indeed, the connected components of an open subset of $\ell^p(\mathbb N)$ are open, since each point has a convex ball contained in that open set. If $D$ is a component of $\{u\mid b(u)\neq0\}$, then $b=0$ on $\partial D$: otherwise, a connected ball on which $b$ is nonzero would meet $D$ and lie in the same component, contradicting the fact that its center is a boundary point. Thus the closure of $b(\partial D)$ is $\{0\}$, and the cited principle gives $b(D)=\{0\}$, a contradiction. The result shows that the differentiability order in the preceding corollary is optimal. The assertion about partitions follows from
Lemma~\ref{lem:cortes-modelo-particiones-banach}, since a partition with this
regularity would produce a bump function of the same order.

There are nonreflexive examples whose construction requires more refined
techniques. Haydon proved the following result.

\begin{definition}[The spaces $C(K)$ and $C_0(L)$]
\label{def:espacios-CK-C0L}
If $K$ is compact and Hausdorff, $C(K)$ is the real Banach space of
continuous functions $u\colon K\longrightarrow\mathbb R$, with the norm
$\displaystyle\|u\|_\infty:=\displaystyle\sup_{x\in K}|u(x)|$. If $L$ is locally compact and
Hausdorff, define
\[
C_0(L)
:=\left\{u\colon L\longrightarrow\mathbb R\middle|
u\text{ is continuous and, for every }\varepsilon>0,
\ \{x\in L\mid |u(x)|\geq\varepsilon\}\text{ is compact}
\right\},
\]
with the supremum norm. This norm is finite: for $\varepsilon=1$, we have $|u|<1$ outside
a compact set and, if this set is nonempty, $|u|$ attains a maximum on it. Moreover, $C_0(L)$ is complete. Indeed, a Cauchy sequence
for $\|\cdot\|_\infty$ converges uniformly to a continuous bounded
function $u$. Given $\varepsilon>0$, if
$\|u-u_n\|_\infty<\frac{\varepsilon}{2}$, then
\[
 \{x\in L\mid |u(x)|\geq\varepsilon\}
 \subseteq
 \{x\in L\mid |u_n(x)|\geq\frac{\varepsilon}{2}\};
\];
the set on the left is closed and the one on the right is compact, so
the former is compact. Its elements are precisely the continuous
functions vanishing at infinity.
\end{definition}

\begin{theorem}[Haydon]
\label{teo:haydon-particiones-Ck}
If $K$ is a countable compact Hausdorff space, then $C(K)$ admits a
smooth norm $\nu$ away from the origin and constants $c_K,C_K>0$, which may depend
on $K$, such that
$c_K\|u\|_\infty\leq\nu(u)\leq C_K\|u\|_\infty$ for every $u\in C(K)$;
moreover, it admits smooth partitions of unity. The same
holds for $C_0([0,\Omega))$, where $\Omega$ is an ordinal and $[0,\Omega)$
has the order topology, with constants that may depend on $\Omega$.
\end{theorem}

The smooth norm is constructed in Haydon's Theorem~1 and Proposition~2,
pp.~457--463; the partition criterion is his
Theorem~2, and the two cases in the statement follow from
Corollaries~4--5, pp.~464--468~\cite{Haydon1996}. Haydon's construction uses Talagrand operators, projectional resolutions, and transfinite induction on ordinal spaces. These examples show that the model space assumption in Theorem~\ref{teo:particiones-unidad-variedad-banach} can also hold without reflexivity.

\subsection{Geometric consequences}

Partitions of unity allow constructions with values in vector spaces
to be globalized. We record a form that will be used
frequently.

\begin{proposition}[Gluing functions with linear codomain]
\label{prop:pegado-lineal-particion-banach}
Let $M$ be a $C^k$--paracompact Banach manifold, let
$(U_\alpha)_{\alpha\in A}$ be an open cover, and let $F$ be a Banach
space. Suppose $f_\alpha\colon U_\alpha\longrightarrow F$ is of class $C^k$ for each
$\alpha$. If $(\psi_\alpha)_{\alpha\in A}$ is a partition of unity of
class $C^k$ subordinate to the cover, then
\[
f\colon M\longrightarrow F,
\qquad
f(p):=\sum_{\alpha\in A}\psi_\alpha(p)f_\alpha(p),
\],
where each product is extended by zero outside $U_\alpha$, is well
defined and of class $C^k$.
\end{proposition}

\begin{proof}
Inclusion $\supp(\psi_\alpha)\subseteq U_\alpha$ ensures that
$\psi_\alpha f_\alpha$ extends by zero to a function of class $C^k$ on
all of $M$. Each point has a neighborhood meeting only finitely
many supports $\supp(\psi_\alpha)$. The sum defining $f$ is therefore
locally finite and locally agrees with a finite sum of functions of
class $C^k$.
\end{proof}

\begin{proposition}[Fine smooth approximation]
\label{prop:aproximacion-suave-fina-banach}
Let $M$ be a $C^k$--paracompact Banach manifold, let $F$ be a Banach
space, and let $f\colon M\longrightarrow F$ be continuous. For every continuous function
$\varepsilon\colon M\longrightarrow(0,\infty)$, there exists $g\in C^k(M,F)$ such that
\[
\|g(p)-f(p)\|_F<\varepsilon(p)
\qquad
\text{for every }p\in M.
\]
\end{proposition}

\begin{proof}
For each $x\in M$, continuity of $f$ and $\varepsilon$ allows us to choose
an open neighborhood $U_x$ of $x$ such that, for every $p\in U_x$,
\[
\|f(p)-f(x)\|_F<\frac{\varepsilon(x)}{2},
\qquad
\varepsilon(p)>\frac{\varepsilon(x)}{2}.
\]
Let $(\psi_x)_{x\in M}$ be a partition of unity of class $C^k$
subordinate to $(U_x)_{x\in M}$. Define
\[
g(p):=\sum_{x\in M}\psi_x(p)f(x).
\]
The sum is locally finite, and each summand is the product of a function of
class $C^k$ and a fixed vector in $F$; thus $g$ is of class $C^k$. If
$\psi_x(p)\neq0$, then $p\in\supp(\psi_x)\subseteq U_x$, and the
preceding inequalities give
\[
\|f(x)-f(p)\|_F
<
\frac{\varepsilon(x)}{2}
<
\varepsilon(p).
\]
Using $\psi_x(p)\geq0$ and $\displaystyle\sum_{x\in M}\psi_x(p)=1$, we conclude that
\[
\|g(p)-f(p)\|_F
\leq
\sum_{x\in M}\psi_x(p)\|f(x)-f(p)\|_F
<
\varepsilon(p).
\]
\end{proof}

In particular, partitions of class $C^k$ imply fine approximation
of continuous functions with linear codomain. Together with
Proposition~\ref{prop:norma-suave-produce-flan-banach} and
Theorem~\ref{teo:criterio-fry-particiones-banach}, this makes precise the
implications used here: if a norm that is $C^k$ away from the origin satisfies
$c\|x\|\leq\nu(x)\leq C\|x\|$ with $c,C>0$, it produces a bump
function; under the Asplund assumption on the dual, this function
produces $C^k$ partitions; and partitions produce smooth approximations.
The converse implications are not needed; they belong to substantially more delicate
renorming and approximation problems.

The same argument preserves linear affine relations whenever the
local objects take values in a common Banach space. Indeed, if
$L\colon F\longrightarrow G$ is continuous and linear, $b\in G$, and the local functions
$f_\alpha\colon U_\alpha\longrightarrow F$ satisfy
$L\circ f_\alpha=b$, then the glued function
$\displaystyle f:=\displaystyle\sum_{\alpha\in A}\psi_\alpha f_\alpha$ satisfies
\[
L\circ f
=
\sum_{\alpha\in A}\psi_\alpha L(f_\alpha)
=
\left(\sum_{\alpha\in A}\psi_\alpha\right)b
=b.
\]
Thus partitions of unity allow local constructions with linear codomain
to be globalized without destroying the affine constraints they
satisfy.

There is no analogous operation for maps with values in an arbitrary Banach
manifold, since a convex combination of points of a manifold is not
defined. Globalizing such maps therefore requires
additional geometric constructions. This distinction appears in the
global study of submersions and mapping spaces; a modern
application of smooth partitions in this setting is the
Stacey--Roberts lemma for Banach manifolds, presented in
\cite{KristelSchmeding2025}.

\section{Fréchet manifolds}
\label{sec:variedades-frechet-autocontenidas}

Constructions by charts and curves require only the chain rule and
continuity of linear differentials. They therefore survive the passage
from Banach to Fréchet spaces. The same is not true of local normal forms
derived from the inverse function theorem. We separate these issues from
the outset.

\begin{definition}[Fréchet atlas]

\label{def:variedad-frechet}

A \textbf{smooth Fréchet manifold without boundary} is a Hausdorff
space $\mathcal M$ equipped with a maximal atlas of homeomorphisms
$\varphi_\alpha\colon U_\alpha\to V_\alpha$, where $V_\alpha$ is open
in a Fréchet space $E_\alpha$, the $U_\alpha$ cover $\mathcal M$, and
all coordinate changes and their inverses are $C^\infty_{\mathrm B}$.
A map between such manifolds is smooth when its chart representations
are $C^\infty_{\mathrm B}$. Atlases of class $C^r_{\mathrm B}$ can also
be defined, but when the model is Banach and $r$ is finite this notation
does not replace the normed class $C^r$ in
Definition~\ref{def:variedad-banach}.

As in our definition of a Banach manifold, neither second countability
nor paracompactness is included among the local axioms. They will be
assumed explicitly when needed. A diffeomorphism is a smooth bijection
with smooth inverse.

\end{definition}

\begin{proposition}[Independence of charts and products]

\label{prop:atlas-productos-frechet}

Every compatible smooth atlas has a unique maximal extension.
Smoothness of a map is independent of the chosen charts, and compositions
of smooth maps are smooth. Open subsets of Fréchet spaces and finite
products of Fréchet manifolds are Fréchet manifolds. Every smooth Banach
manifold is a smooth Fréchet manifold with the same charts.

\end{proposition}

\begin{proof}

Adjoin every chart compatible with the given atlas. If two adjoined
charts intersect, choose an original chart near any point of their
intersection. The transition between the adjoined charts is a composition
of two compatible transitions through the original chart. It is smooth
by Theorem~\ref{teo:reglas-calculo-frechet}, as is its inverse.
This proves compatibility, maximality, and uniqueness. Independence of
chart representations and the composition rule follow from the same
factorization.

For a product $E\times F$, the seminorms $p_j(x)+q_j(y)$ make it
Fréchet: a Cauchy sequence converges in each factor. Products of charts
have product coordinate changes, which are smooth by definition.
The final assertion follows from agreement of the smooth classes proved
in Proposition~\ref{prop:comparacion-calculos-banach-frechet}.

\end{proof}

\subsection{Tangents, differentials, and bundles}
\label{subsec:tangentes-haces-frechet}

\begin{theorem}[The tangent bundle of a Fréchet manifold]

\label{teo:tangente-frechet}

Let $\mathcal M$ be a smooth Fréchet manifold. Define $T_p\mathcal M$
by $C^1_{\mathrm B}$ curves through $p$: two curves are equivalent when
their velocities agree in one chart. This definition is independent of
the chart. Every chart $(U,\varphi)$ induces a linear topological
isomorphism

\[
 \Theta_{\varphi,p}\colon T_p\mathcal M\longrightarrow E,
 \qquad[\gamma]\longmapsto(\varphi\circ\gamma)'(0).
\]

The disjoint union $T\mathcal M$ is a Fréchet manifold with charts
$(p,v)\mapsto(\varphi(p),\Theta_{\varphi,p}v)$ and smooth projection
$\pi\colon T\mathcal M\to\mathcal M$.
If $F\colon\mathcal M\to\mathcal N$ is smooth, the map
$dF_p[\gamma]=[F\circ\gamma]$ is well defined, linear, and continuous,
and the differentials form a smooth map
$TF\colon T\mathcal M\to T\mathcal N$. One has $T(G\circ F)=TG\circ TF$.

\end{theorem}

\begin{proof}

The chain rule shows that, for $h=\psi\circ\varphi^{-1}$, velocities
transform by $Dh(\varphi(p))$. This operator is a continuous linear
isomorphism with inverse $D(h^{-1})(\psi(p))$, since differentiating the
inverse identities gives the linear identities. Thus the relation between
curves is independent of the chart and is an equivalence relation.

Every $v\in E$ is represented by
$t\mapsto\varphi^{-1}(\varphi(p)+tv)$ for small $t$, proving
surjectivity of $\Theta_{\varphi,p}$; injectivity is the definition of
equivalence. Transport the operations and topology of $E$.
The preceding linear topological coordinate changes show that these
structures are independent of the chart. On the total space, transitions
are

\[
 (x,v)\longmapsto\bigl(h(x),Dh(x)v\bigr).
\]

They are smooth by the derivative evaluation assertion in
Theorem~\ref{teo:reglas-calculo-frechet}. They give a Hausdorff topology:
vectors over different base points are separated using the projection,
whereas vectors over the same base point are separated in a trivialization.
The projection is locally $(x,v)\mapsto x$ and hence smooth.
Finally, in charts $TF$ has the same form, with $h$ replaced by the chart
representation of $F$. The chain rule proves the final identity.

\end{proof}

\begin{definition}[Fréchet vector bundles]

\label{def:haz-vectorial-frechet}

A \textbf{smooth Fréchet vector bundle} over $\mathcal M$ is a map
$\pi\colon\mathbf E\to\mathcal M$ with trivializations
$\tau_\alpha\colon\pi^{-1}(U_\alpha)\to U_\alpha\times F_
\alpha$ that are diffeomorphisms and linear on each fiber. Transitions
have the form $(p,v)\mapsto(p,g_{\beta\alpha}(p)v)$; this map and its
inverse are required to be jointly smooth. Here $g_{\beta\alpha}(p)$ is
a linear topological isomorphism. This formulation does not assume that
the set of isomorphisms is open in an operator space, a property not
available in general for Fréchet spaces.

The preceding theorem constructs $T\mathcal M$ as such a bundle.
A smooth section of $T\mathcal M$ is called a vector field.
The definition does not assert that every vector field has a flow.

\end{definition}

\begin{remark}[The dual is not another Fréchet model by definition]

\label{obs:dual-frechet-no-automatico}

In the chapter on Banach bundles, $T'\mathcal M$ is constructed using
the continuous dual of the model with its operator norm. We do not transfer
that construction unchanged to Fréchet spaces. For example, a continuous
linear functional on $E=\mathbb R^{\mathbb N_0}$ depends on finitely many
coordinates: the linear continuity criterion gives
$|\ell(x)|\leq Cp_j(x)$, so $\ell$ vanishes when the first $j+1$
coordinates vanish.

Moreover, no vector-space topology on this dual makes evaluation
$E\times E'\to\mathbb R$ jointly continuous. Indeed, continuity would
give zero neighborhoods $U$ and $V$ such that $|\ell(x)|<1$ for
$x\in U$ and $\ell\in V$. The set $U$ contains a basic neighborhood
restricting only the first $j+1$ coordinates. Choose $n>j$.
Since $V$ is absorbing, it contains $a e_n'$ for some $a>0$, where
$e_n'(x)=x_n$. Yet $te_n\in U$ for every $t$, and
$ae_n'(te_n)=at$ is not bounded by $1$.
We therefore define forms through their smooth scalar evaluations,
not through an assumed universal Fréchet cotangent bundle.

\end{remark}

\subsection{Submanifolds and a theorem for finite constraints}
\label{subsec:restricciones-finitas-frechet}

\begin{definition}[Split submanifold]

\label{def:subvariedad-frechet-complementada}

A subset $\mathcal S\subseteq\mathcal M$ is a \textbf{split Fréchet
submanifold} if every one of its points has a chart modeled on
$E=F\oplus G$, a topological direct sum of closed subspaces, carrying
$\mathcal S$ to the intersection of the coordinate open set with
$F\times\{0\}$. Closed subspaces of a Fréchet space are complete for
the restricted seminorms. The induced structure on $\mathcal S$ and its
tangent space are obtained by restricting these charts, as in
Proposition~\ref{prop:subvariedad-hereda-estructura-banach}.

In particular, the graph of a smooth map $h\colon U\subseteq F\to G$
is a split submanifold: the diffeomorphism
$(x,y)\mapsto(x,y-h(x))$, with inverse $(x,z)\mapsto(x,z+h(x))$,
carries it to $U\times\{0\}$. This assertion has an explicit proof and
does not use an abstract Fréchet inverse theorem.

\end{definition}

\begin{theorem}[Finite-codimensional level sets]

\label{teo:nivel-finito-frechet}

Let $f\colon U\subseteq E\to\mathbb R^q$ be smooth, with $E$ Fréchet.
If $Df(x_0)$ is surjective, then $f^{-1}(f(x_0))$ is, near $x_0$, a split
submanifold with tangent space $\ker Df(x_0)$.
The same conclusion holds in charts on a Fréchet manifold.

\end{theorem}

\begin{proof}

After translation, assume $x_0=0$ and $f(0)=0$. Choose
$v_1,\ldots,v_q\in E$ mapped by $Df(0)$ to the standard basis.
The operator $R(a)=\displaystyle \sum_{i=1}^q a_i v_i$ is continuous and satisfies
$Df(0)R=I$. Thus $E=K\oplus R(\mathbb R^q)$, with
$K=\ker Df(0)$, and we write $F(z,y)=f(z+Ry)$.
Then $D_yF(0,0)=I$.

Joint continuity of $D_yF$ implies matrix-norm continuity because
its direction space is finite-dimensional. There are a zero neighborhood
$V$ in $K$ and $\rho>0$ such that
$\|I-D_yF(z,y)\|\leq1/2$ for $z\in V$ and $|y|\leq\rho$.
Shrinking $V$ also gives $|F(z,0)|<\rho/4$.
The map $T_z(y)=y-F(z,y)$ maps the closed ball of radius $\rho$ into
itself and is a contraction with constant $1/2$. Its unique fixed point
$y=h(z)$ satisfies $|h(z)|<\rho/2$ and $F(z,h(z))=0$.
All solutions in that ball have this form.

The contraction estimate gives

\[
 |h(z')-h(z)|\leq2|F(z',h(z))-F(z,h(z))|,
\]

and proves continuity of $h$. For $z'=z+tv$, put
$\Delta h=h(z+tv)-h(z)$. The integral formula, first in $z$ and then
in $y$, gives

\[
 0=t\int_0^1D_zF(z+stv,h(z+tv))v\,ds
   +\left(\int_0^1D_yF(z,h(z)+s\Delta h)\,ds\right)\Delta h.
\]

The matrix in parentheses is invertible, with inverse norm at most
$2$, and converges to $D_yF(z,h(z))$. Dividing by $t$ gives

\[
 Dh(z)v=-D_yF(z,h(z))^{-1}D_zF(z,h(z))v.
\]

The right-hand side is jointly continuous in $(z,v)$, so $h$ is
$C^1_{\mathrm B}$. Matrix inversion is smooth; the chain rule and
derivative evaluation maps prove inductively that $h$ is
$C^r_{\mathrm B}$ for every $r$. The level set is the smooth graph
$y=h(z)$. Since $D_zF(0,0)=0$, one has $Dh(0)=0$, so the tangent
space to that graph is $K$. The construction only inverted matrices and
contracted a finite-dimensional ball; it does not prove an analogous
theorem for an arbitrary Fréchet target.

\end{proof}

\subsection{Differential forms without a preassigned dual bundle}
\label{subsec:formas-frechet}

\begin{definition}[Smooth differential form]

\label{def:forma-diferencial-frechet}

A smooth $k$-form $\omega$ on $\mathcal M$ assigns to each $p$ a
continuous alternating $k$-linear map
$\omega_p\colon(T_p\mathcal M)^k\to\mathbb R$ such that, in each chart,
the evaluation map $(x,v_1,\ldots,v_k)\mapsto
\omega_x(v_1,\ldots,v_k)$ is jointly smooth. For $k=0$, this means a
smooth function. Denote the space of these forms by $\Omega^k(\mathcal M)$.
Pullback by $F\colon\mathcal M\to\mathcal N$ is defined by

\[
 (F^*\omega)_p(v_1,\ldots,v_k)
 =\omega_{F(p)}(dF_pv_1,\ldots,dF_pv_k).
\]

We use the same wedge-product convention as in finite dimensions:

\[
 (\omega\wedge\eta)_p(v_1,\ldots,v_{k+\ell})
 =\frac1{k!\ell!}\displaystyle\sum_{\sigma\in S_{k+\ell}}
 \operatorname{sign}(\sigma)
 \omega_p(v_{\sigma(1)},\ldots,v_{\sigma(k)})
 \eta_p(v_{\sigma(k+1)},\ldots,v_{\sigma(k+\ell)}).
\]

\end{definition}

\begin{theorem}[Exterior derivative and Poincaré lemma]

\label{teo:exterior-poincare-frechet}

There is an operator $d\colon\Omega^k(\mathcal M)\to
\Omega^{k+1}(\mathcal M)$ given in charts by

\begin{equation}\label{eq:exterior-frechet-cartas}
 (d\omega)_x(v_0,\ldots,v_k)
 =\displaystyle\sum_{i=0}^k(-1)^i
 D_x\bigl(\omega_{\,\cdot\,}(v_0,\ldots,\widehat v_i,\ldots,v_k)\bigr)(x)v_i.
\end{equation}

It is independent of the charts and satisfies $d^2=0$,
$d(\omega\wedge\eta)=d\omega\wedge\eta+(-1)^k\omega\wedge d\eta$,
and $dF^*=F^*d$. If $U\subseteq E$ is open and star-shaped about $0$,
every closed form of degree $k\geq1$ on $U$ is exact.

\end{theorem}

\begin{proof}

Expression \eqref{eq:exterior-frechet-cartas} is alternating and its
evaluation is smooth by the calculus of derivative evaluation maps.
In $d(d\omega)$, each pair of derivatives in directions $v_i,v_j$
appears twice with opposite signs. Symmetry of second derivatives makes
them cancel, proving $d^2=0$.

For naturality, differentiate
$\omega_{F(x)}(DF(x)v_1,\ldots,DF(x)v_k)$.
Terms differentiating $\omega$ give $F^*(d\omega)$ by the chain rule.
The remaining terms contain $D^2F(x)(v_i,v_j)$; the alternating sum
contains another term with $i,j$ interchanged and opposite sign.
They cancel by symmetry of $D^2F$. Hence $dF^*=F^*d$ on model open
sets. Applying this to coordinate changes glues the local definitions.
The wedge product can be written as a sum over shuffles by grouping the
$k!\ell!$ internal permutations. On differentiating it, separate terms
where the derivative hits the first or second factor. The former give
$d\omega\wedge\eta$; in the latter, the differentiating direction passes
the $k$ entries of the first factor, producing $(-1)^k$.
This proves the Leibniz rule with the stated convention.

On the star-shaped open set, define, for $k\geq1$,

\[
 (K\omega)_x(v_1,\ldots,v_{k-1})
 =\int_0^1t^{k-1}\omega_{tx}(x,v_1,\ldots,v_{k-1})\,dt.
\]

We make smoothness of this evaluation explicit. Fix $x\in U$ and
choose an open convex neighborhood $V$ of $x$ contained in $U$.
Since $U$ is star-shaped, $ty\in U$ for every $(t,y)\in[0,1]\times V$.
On this domain, define
\[
 F(t,z):=t^{k-1}\omega_{ty}(y,w_1,\ldots,w_{k-1}),\qquad
 z=(y,w_1,\ldots,w_{k-1})\in V\times E^{k-1}.
\]
For fixed $z$ and a direction $\eta\in E^k$, the integral mean-value
formula gives, for sufficiently small $s\neq0$,
\[
 \frac{F(t,z+s\eta)-F(t,z)}{s}
 =\int_0^1D_zF(t,z+as\eta)\eta\,da.
\]
The integrand converges to $D_zF(t,z)\eta$ uniformly in
$(t,a)\in[0,1]^2$. To check this, given $\varepsilon>0$, joint continuity
of $D_zF$ gives, around each $(t,z,\eta)$, an interval $I_t$ and
neighborhoods of $z$ and $\eta$ such that
\[
 |D_zF(\tau,z')\eta'-D_zF(\tau,z)\eta|<\varepsilon
 \quad\text{if }\tau\in I_t
\]
for $z',\eta'$ in these neighborhoods. A finite subcover of $[0,1]$ and the intersection
of these neighborhoods give the same bound for every $t$.
Since $z+as\eta\to z$ uniformly in $a$, it follows that
\[
 \sup_{t,a\in[0,1]}
 |D_zF(t,z+as\eta)\eta-D_zF(t,z)\eta|\longrightarrow0.
\]
We may therefore integrate the difference quotient. Applying the same
argument to every $D_z^jF$ and its direction variables proves inductively
\[
 D^j\left(\int_0^1F(t,\cdot)\,dt\right)(z)
       (\eta_1,\ldots,\eta_j)
 =\int_0^1D_z^jF(t,z)(\eta_1,\ldots,\eta_j)\,dt,
\]
with joint continuity in all variables. This proves smoothness of
$K\omega$ in the Bastiani sense.

To compute the homotopy identity, write
$(D\omega)_y(a;b_1,\ldots,b_k)$ for the derivative, in the direction $a$,
of $y\mapsto\omega_y(b_1,\ldots,b_k)$, and set
\[
 B_i(t):=(D\omega)_{tx}
       (v_i;x,v_1,\ldots,\widehat v_i,\ldots,v_k).
\]
Formula \eqref{eq:exterior-frechet-cartas} and differentiation under
the integral give
\begin{align*}
 (dK\omega)_x(v_1,\ldots,v_k)
 &=\int_0^1\left[
 k t^{k-1}\omega_{tx}(v_1,\ldots,v_k)
 +t^k\sum_{i=1}^k(-1)^{i-1}B_i(t)\right]dt,\\
 (Kd\omega)_x(v_1,\ldots,v_k)
 &=\int_0^1t^k\left[
 (D\omega)_{tx}(x;v_1,\ldots,v_k)
 +\sum_{i=1}^k(-1)^iB_i(t)\right]dt.
\end{align*}
The terms $B_i$ cancel upon addition. The remaining integrand is
\[
 k t^{k-1}\omega_{tx}(v_1,\ldots,v_k)
 +t^k(D\omega)_{tx}(x;v_1,\ldots,v_k)
 =\frac{d}{dt}\bigl[t^k\omega_{tx}(v_1,\ldots,v_k)\bigr].
\]
For $k\geq1$, the boundary term at $t=0$ is zero. Consequently,

\[
 dK\omega+Kd\omega=\omega\quad(k\geq1),
 \qquad Kdf=f-f(0)\quad(k=0).
\]

The second identity is the fundamental theorem of calculus.
If $d\omega=0$ and $k\geq1$, the first gives $\omega=d(K\omega)$.

\end{proof}

\begin{corollary}[Homotopy and de Rham cohomology]

\label{cor:homotopia-derham-frechet}

The quotients
$H^k_{\mathrm{dR}}(\mathcal M)=\ker d\vert_{\Omega^k}/
\operatorname{Ran}d\vert_{\Omega^{k-1}}$ are defined, and smoothly
homotopic maps induce the same map on them.

\end{corollary}

\begin{proof}

Let $H\colon[0,1]\times\mathcal M\to\mathcal N$ be a homotopy
smooth up to the endpoints. On the product, write
$H^*\omega=dt\wedge\alpha_t+\beta_t$, where
$\beta_t=H_t^*\omega$ and $\alpha_t$ is obtained by inserting
$\partial_t$ in the first entry. The local formula for $d$ gives

\[
 d(H^*\omega)=dt\wedge(\partial_t\beta_t-d\alpha_t)+d\beta_t.
\]

Define $K_H\omega=\displaystyle \int_0^1\alpha_t\,dt$ and integrate this identity
to obtain
$H_1^*\omega-H_0^*\omega=dK_H\omega+K_Hd\omega$.
Smoothness and differentiation under the integral are justified as in
the preceding theorem. For closed forms, the difference is exact.

\end{proof}

This construction gives the de Rham complex and invariance under
smooth homotopies. By itself, it does not identify these groups with
singular or Čech cohomology on an arbitrary Fréchet manifold.
Such an identification requires global hypotheses and an additional proof.
Nor have smooth partitions of unity been attributed to all Fréchet models.
The partition hypotheses already specified in the Banach category retain
their role in global results.

\chapter{Vector bundles, connections, and flows on Banach manifolds}
\label{cap:haces-flujos-banach}

Vector bundles over Banach manifolds are described by transition
maps with values in groups of invertible operators. The
topology used on these groups is induced by the operator norm;
this choice makes Fréchet calculus applicable and avoids confusing the
joint continuity of evaluation with the regularity, in the operator
norm, of a family of linear transformations.

The framework of bundles with Banach fibers and the need to distinguish
topological tensor products already appear in the classical formulation of
infinite-dimensional global analysis
\cite[\S\S4\textup{(C)}--4\textup{(G)}]{Eells1966}. Here we work in the
Banach category and use Fréchet differentiability.

\begin{semblanzaHistorica}{Differential geometry with infinite-dimensional fibers}
Banach bundles allow each fiber to be a complete infinite-dimensional space. This generalization arises naturally in the study of families of sections, operators, or initial data that vary over a base. The operator norm determines the regularity of transition functions, connections, and parallel transport. We will study these constructions together with the topological conditions they require, some of which hold automatically in finite dimensions.
\end{semblanzaHistorica}

\section{Banach vector bundles}
\label{sec:haces-vectoriales-banach}

Let $X$ and $F$ be real Banach spaces, and let $M$ be a Banach manifold of
class $C^r$, with $r\in\mathbb N\cup\{\infty\}$, modeled on $X$. The
space $F$ will be the model fiber. We use the convention of
Definition~\ref{def:variedad-banach}: charts have open images in
the model spaces. A continuous bundle is understood to satisfy
the same conditions as a $C^r$ bundle, with $C^r$ regularity
replaced by continuity in the operator norm. This variant will be used for the
tangent bundle of a $C^1$ manifold.

\begin{proposition}[General linear group of a Banach space]
\label{prop:grupo-lineal-general-fibra-banach}
The set
\[
\operatorname{GL}(F)
:=
\{A\in\mathcal L(F)\mid A\text{ is invertible}\}
\]
is an open subset of the Banach algebra $\mathcal L(F)$. With this structure,
$\operatorname{GL}(F)$ is a Banach Lie group: composition and
inversion
\[
(A,B)\longmapsto AB,
\qquad
A\longmapsto A^{-1}
\]
are smooth. Moreover,
\[
D(A\mapsto A^{-1})(A)H=-A^{-1}HA^{-1}.
\]
\end{proposition}

\begin{proof}
Composition is continuous and bilinear, and hence smooth. The openness of the
set of invertible operators, the smoothness of inversion, and the formula
for its differential are precisely the assertions of
Proposition~\ref{prop:inversion-operadores-suave}, applied with $X=Y=F$.
\end{proof}

Before taking duals of fibers, it is useful to distinguish two notions that coincide in
finite dimensions but not in a general Banach space. For a vector
space $V$, we denote its algebraic dual, consisting of all
linear functionals, by $V^*$. For a Banach space $F$, we denote its
continuous dual, consisting of the continuous linear functionals, by
\[
F':=\mathcal L(F;\mathbb R)
\]
The latter is a subspace of the algebraic dual $F^*$ and carries the norm
\[
\|\lambda\|_{F'}
:=
\sup_{v\in F,\ \|v\|_F\leq1}|\lambda(v)|.
\]
\index{dual!algebraic}\index{dual!continuous}

\begin{proposition}[Continuous dual and transpose operators]
\label{prop:dual-topologico-banach}\glsadd{dual-topologico-banach}
The space $F'$ is a Banach space. If $A\in\mathcal L(F,G)$, its
transpose
\[
A'\colon G'\longrightarrow F',
\qquad
A'\mu:=\mu\circ A,
\]
is linear and continuous, and satisfies
\[
\|A'\|_{\mathcal L(G',F')}=\|A\|_{\mathcal L(F,G)}.
\]
If $A\colon F\longrightarrow G$ is a topological linear isomorphism, then
so is $A'\colon G'\longrightarrow F'$, and $(A')^{-1}=(A^{-1})'$. Moreover, the map
\[
\mathcal L(F,G)\longrightarrow\mathcal L(G',F'),
\qquad
A\longmapsto A',
\]
is an isometric linear embedding.

If $F$ is infinite-dimensional, the inclusion $F'\subseteq F^*$ is
strict.
\end{proposition}

\begin{proof}
Let $(\lambda_n)_{n\in\mathbb N}$ be a Cauchy sequence in $F'$. For each $v\in F$, the
sequence $(\lambda_n(v))_{n\in\mathbb N}$ is Cauchy in $\mathbb R$, since
\[
|\lambda_n(v)-\lambda_m(v)|
\leq
\|\lambda_n-\lambda_m\|_{F'}\|v\|_F.
\]
Define $\displaystyle\lambda(v):=\lim_{n\to\infty}\lambda_n(v)$. Passing to the limit shows that
$\lambda$ is linear. Since every Cauchy sequence is bounded, there exists $C>0$
such that $\|\lambda_n\|_{F'}\leq C$ for every $n\in\mathbb N$, and then
$|\lambda(v)|\leq C\|v\|_F$. Therefore, $\lambda\in F'$. Given
$\varepsilon>0$, for sufficiently large $n,m$ we have
$\|\lambda_n-\lambda_m\|_{F'}<\varepsilon$; letting $m$ tend to infinity
in the corresponding pointwise inequality gives
$\|\lambda_n-\lambda\|_{F'}\leq\varepsilon$. Thus, $F'$ is complete.

The inequality $\|A'\|\leq\|A\|$ follows from
\[
|(A'\mu)(v)|
\leq
\|\mu\|_{G'}\|A\|\|v\|_F.
\]
Conversely, the Hahn--Banach theorem applied to $Av\in G$ yields, for
each $v\in F$ with $Av\neq0$, a functional $\mu_v\in G'$ such that
$\|\mu_v\|_{G'}=1$ and $\mu_v(Av)=\|Av\|_G$. Hence,
\[
\|Av\|_G
\leq
\|A'\|\|v\|_F,
\]
an inequality that also holds when $Av=0$. Taking the supremum over the
unit ball gives $\|A\|\leq\|A'\|$. The assertions about the inverse
follow from $(BA)'=A'B'$.

Finally, suppose that $F$ is infinite-dimensional. There exists a linearly
independent sequence $(v_n)_{n\in\mathbb N}$ of unit vectors. On its linear
span, consisting of finite linear combinations of these vectors,
we define $\lambda_0(v_n):=n$ for each $n\in\mathbb N$ and extend it
algebraically to $F$ by extending $(v_n)_{n\in\mathbb N}$ to a Hamel basis. The extension
$\lambda\in F^*$ is not continuous, because it is not bounded on the unit
sphere. Consequently, $\lambda\notin F'$.
\end{proof}

The algebraic dual does not in general carry a canonical Banach topology
compatible with evaluation. For this reason, all the dual bundles that
appear below are constructed using the continuous dual.

\begin{definition}[Banach vector bundle]
\label{def:haz-vectorial-banach}
Let $M$ be a Banach manifold of class $C^r$.
\index{vector bundle!Banach}
A \textbf{Banach vector bundle of class $C^r$ with model fiber $F$} over
$M$ consists of a Hausdorff topological
space $\boldsymbol{\mathcal{E}}$, a continuous surjective map
\[
\pi\colon \boldsymbol{\mathcal{E}}\longrightarrow M
\]
and a Banach vector space structure on each fiber
$\boldsymbol{\mathcal{E}}_p:=\pi^{-1}(p)$, together with an open cover
$(U_\alpha)_{\alpha\in A}$ of $M$ and homeomorphisms
\[
\Phi_\alpha\colon \pi^{-1}(U_\alpha)\longrightarrow U_\alpha\times F
\]
satisfying the following conditions:
\begin{enumerate}[label=(\alph*)]
\item $\operatorname{pr}_1\circ\Phi_\alpha=\pi$;
\item for each $p\in U_\alpha$, the restriction
\[
\Phi_{\alpha,p}
:=
\operatorname{pr}_2\circ
\bigl(\Phi_\alpha\restriction_{\boldsymbol{\mathcal{E}}_p}\bigr)
\colon \boldsymbol{\mathcal{E}}_p\longrightarrow F
\]
is a linear and topological isomorphism;
\item on every nonempty intersection there are maps of class $C^r$, with respect to
the operator norm topology,
\[
g_{\beta\alpha}\colon U_\alpha\cap U_\beta
\longrightarrow\operatorname{GL}(F)
\]
such that
\[
\Phi_\beta\circ\Phi_\alpha^{-1}(p,v)
=
\bigl(p,g_{\beta\alpha}(p)v\bigr).
\]
\end{enumerate}
The family $(U_\alpha,\Phi_\alpha)$ is called a
\textbf{family of local trivializations}. Two families determine the same
$C^r$ bundle if their union satisfies condition \textup{(c)}.
\end{definition}

The regularity in Definition~\ref{def:haz-vectorial-banach} is required in the
operator norm. Since evaluation
\[
\mathcal L(F)\times F\longrightarrow F,
\qquad
(A,v)\longmapsto Av,
\]
is continuous and bilinear, each change of trivialization is of class $C^r$ in the
base and fiber coordinates. Combined with an atlas of $M$, the
trivializations equip $\boldsymbol{\mathcal E}$ with a Banach manifold atlas
modeled on $X\times F$. The projection $\pi$ is then of class $C^r$.

Definition~\ref{def:haz-vectorial-banach} reproduces, with Banach fibers
and regularity in the operator norm, the finite-dimensional construction in
Definition~\ref{def:nociones-fundamentales-haz-vectorial-suave-de-rango-sobre}.
In finite rank, continuity of a family of matrices and continuity in the
operator norm are the same condition; for a Banach fiber, the
latter is the condition that allows Fréchet calculus to be applied directly.

\begin{proposition}[Transition maps and reconstruction]
\label{prop:transiciones-haz-banach}
Let $M$ be a Banach manifold of class $C^r$.
The transition maps of a Banach vector bundle satisfy
\[
g_{\alpha\alpha}(p)=\operatorname{id}_F,
\qquad
g_{\alpha\beta}(p)=g_{\beta\alpha}(p)^{-1},
\qquad
g_{\gamma\beta}(p)g_{\beta\alpha}(p)=g_{\gamma\alpha}(p)
\]
for $\alpha,\beta,\gamma\in A$ and $p$ in the intersection of the corresponding
domains.

Conversely, an open cover $(U_\alpha)_{\alpha\in A}$ of $M$ and a family of
$C^r$ maps
\[
g_{\beta\alpha}\colon U_\alpha\cap U_\beta\longrightarrow\operatorname{GL}(F)
\]
satisfying these identities determine, up to a bundle isomorphism of class
$C^r$, a Banach vector bundle with model fiber $F$.
\end{proposition}

\begin{proof}
The identities follow by composing the changes of trivialization. For the
converse, consider the disjoint union
\[
\boldsymbol{\mathcal{E}}_0:=\coprod_{\alpha\in A}(U_\alpha\times F)
\]
and the relation
\[
(\alpha,p,v)\sim(\beta,p,w)
\quad\Longleftrightarrow\quad
w=g_{\beta\alpha}(p)v.
\]
The cocycle identities show that this is an equivalence relation.
Let $q\colon \boldsymbol{\mathcal{E}}_0\longrightarrow\boldsymbol{\mathcal{E}}:=\boldsymbol{\mathcal{E}}_0/{\sim}$ be the projection, and
equip $\boldsymbol{\mathcal{E}}$ with the quotient topology. For each $\alpha$, define
\[
j_\alpha\colon U_\alpha\times F\longrightarrow\boldsymbol{\mathcal{E}},
\qquad
j_\alpha(p,v):=q(\alpha,p,v).
\]
The map $j_\alpha$ is injective. Moreover, if
$O\subseteq U_\alpha\times F$ is open, then
\[
q^{-1}\bigl(j_\alpha(O)\bigr)
=
\coprod_{\beta\in A}
\left\{
(\beta,p,g_{\beta\alpha}(p)v)
\ \middle|\
(p,v)\in O,\ p\in U_\alpha\cap U_\beta
\right\}.
\]
For each $\beta$, the set on the right is the image of the open set
$O\cap((U_\alpha\cap U_\beta)\times F)$ under the homeomorphism
\[
(p,v)\longmapsto(p,g_{\beta\alpha}(p)v).
\]
Consequently, $q^{-1}(j_\alpha(O))$ is open, and the definition of the
quotient topology implies that $j_\alpha(O)$ is open. Therefore,
$j_\alpha$ is a homeomorphism onto the open set
$j_\alpha(U_\alpha\times F)$.

The formula
\[
\pi\bigl(j_\alpha(p,v)\bigr):=p
\]
is well defined, is surjective, and is continuous on each of the open sets
$j_\alpha(U_\alpha\times F)$. In fact,
$j_\alpha(U_\alpha\times F)=\pi^{-1}(U_\alpha)$: if
$q(\beta,p,w)$ projects to $p\in U_\alpha$, it is also
represented by $(\alpha,p,g_{\alpha\beta}(p)w)$. We show that
$\boldsymbol{\mathcal{E}}$ is Hausdorff. If
$e_1,e_2\in\boldsymbol{\mathcal{E}}$ project to distinct points, choose disjoint open subsets
of $M$ containing their projections; their preimages under $\pi$
separate $e_1$ and $e_2$. If $\pi(e_1)=\pi(e_2)=p$, choose
$\alpha$ with $p\in U_\alpha$ and write
$e_i=j_\alpha(p,v_i)$. The inequality $e_1\neq e_2$ implies $v_1\neq v_2$.
Two disjoint balls in $F$ about $v_1$ and $v_2$, multiplied by an
open neighborhood of $p$ contained in $U_\alpha$ and transported by
$j_\alpha$, separate the two points.

On each fiber, the operations are defined by transporting those of $F$ via
$j_\alpha$; the cocycle guarantees that they do not depend on $\alpha$. The
trivializations
\[
\Phi_\alpha\colon j_\alpha(U_\alpha\times F)\longrightarrow U_\alpha\times F,
\qquad
\Phi_\alpha:=j_\alpha^{-1},
\]
have as their coordinate changes precisely the maps
$(p,v)\mapsto(p,g_{\beta\alpha}(p)v)$. Combining them with charts of $M$
gives the Banach manifold atlas on the total space.

Finally, if $\widetilde{\boldsymbol{\mathcal{E}}}\longrightarrow M$ is another bundle equipped with
trivializations $\widetilde\Phi_\alpha$ having the same transition
maps, the formula
\[
j_\alpha(p,v)
\longmapsto
\widetilde\Phi_\alpha^{-1}(p,v)
\]
is independent of $\alpha$. It defines a bundle isomorphism of class $C^r$, since
it is the identity in each trivialization. This proves uniqueness up to
isomorphism.
\end{proof}

Two cocycles $(g_{\beta\alpha})$ and
$(\widetilde g_{\beta\alpha})$ over the same cover are called
\textbf{equivalent} if there are maps
$h_\alpha\in C^r(U_\alpha,\operatorname{GL}(F))$ such that
\begin{equation}
\label{eq:cociclos-equivalentes-haces-banach}
\widetilde g_{\beta\alpha}(p)
=
h_\beta(p)g_{\beta\alpha}(p)h_\alpha(p)^{-1}.
\end{equation}
In that case, the formula
\[
[\alpha,p,v]_g
\longmapsto
[\alpha,p,h_\alpha(p)v]_{\widetilde g}
\]
is well defined by \eqref{eq:cociclos-equivalentes-haces-banach} and
provides an isomorphism of the reconstructed bundles. Conversely, the
local representations of a bundle isomorphism over $M$ provide
a family $(h_\alpha)$ satisfying that identity, after passing to a
common cover. Thus, equivalence of cocycles describes exactly the
change of trivializations of the same bundle.

\begin{definition}[Bundle homomorphism]
\label{def:homomorfismo-haces-banach}
Let $M$ and $N$ be Banach manifolds.
Let $\pi_{\boldsymbol{\mathcal{E}}}\colon \boldsymbol{\mathcal{E}}\longrightarrow M$ and
$\pi_{\boldsymbol{\mathcal{H}}}\colon \boldsymbol{\mathcal{H}}\longrightarrow N$ be Banach vector bundles of class
$C^r$, with model fibers $F$ and $G$, respectively. A
\textbf{bundle homomorphism of class $C^r$ over}
$f\in C^r(M,N)$ is a map $\mathscr A\colon \boldsymbol{\mathcal{E}}\longrightarrow\boldsymbol{\mathcal{H}}$ such that
\[
\pi_{\boldsymbol{\mathcal{H}}}\circ\mathscr A=f\circ\pi_{\boldsymbol{\mathcal{E}}},
\]
whose restrictions
\[
\mathscr A_p\colon \boldsymbol{\mathcal{E}}_p\longrightarrow\boldsymbol{\mathcal{H}}_{f(p)}
\]
are linear and continuous and whose local representations are
$C^r$ families, in the operator norm,
\[
A_{\beta\alpha}\colon U_\alpha\cap f^{-1}(V_\beta)
\longrightarrow\mathcal L(F,G).
\]
A homomorphism over $\operatorname{id}_M$ is called a homomorphism over $M$. It is
a \textbf{bundle isomorphism} if every fiber operator is invertible and the
inverse homomorphism is of class $C^r$.
\end{definition}

\begin{proposition}[Local data of a homomorphism]
\label{prop:datos-locales-homomorfismo-haces-banach}
Let $M$ and $N$ be Banach manifolds, and let $\boldsymbol{\mathcal E}\to M$ and $\boldsymbol{\mathcal H}\to N$ be Banach vector bundles of class $C^r$.
The local representations of a bundle homomorphism satisfy
\begin{equation}
\label{eq:compatibilidad-homomorfismo-haces-banach}
A_{\delta\gamma}(p)
=
h_{\delta\beta}(f(p))A_{\beta\alpha}(p)g_{\alpha\gamma}(p),
\end{equation}
whenever all the trivializations involved are defined. Here
$g$ and $h$ are the transition maps of $\boldsymbol{\mathcal E}$ and $\boldsymbol{\mathcal H}$.

Conversely, a family of $C^r$ maps with values in
$\mathcal L(F,G)$ satisfying
\eqref{eq:compatibilidad-homomorfismo-haces-banach} determines a unique bundle
homomorphism of class $C^r$ over $f$. Composition of homomorphisms corresponds to
pointwise composition of their local operators.
\end{proposition}

\begin{proof}
The first identity follows by writing a local representation with respect to
two different pairs of trivializations. For the converse, let
$e\in\boldsymbol{\mathcal E}_p$, let $\Phi_\alpha$ be a trivialization about $p$, and
choose a trivialization $\Psi_\beta$ of $\boldsymbol{\mathcal H}$ about $f(p)$.
If $\Phi_\alpha(e)=(p,v)$, set
\[
\mathscr A(e)
:=
\Psi_\beta^{-1}\bigl(f(p),A_{\beta\alpha}(p)v\bigr).
\]
Identity \eqref{eq:compatibilidad-homomorfismo-haces-banach} shows that the
result does not depend on $\alpha$ or $\beta$. In each pair of
trivializations, $\mathscr A$ is represented by
$(p,v)\mapsto(f(p),A_{\beta\alpha}(p)v)$, which is of class $C^r$ because
evaluation $\mathcal L(F,G)\times F\longrightarrow G$ is continuous and bilinear. Uniqueness and
the assertion about compositions are immediate in these representations.
\end{proof}

\begin{definition}[Complemented subbundle]
\label{def:subhaz-complementado-banach}
Let $M$ be a Banach manifold.
Let $\boldsymbol{\mathcal{E}}\longrightarrow M$ be a Banach vector bundle with model fiber $F$. A
\textbf{complemented vector subbundle} is a subset
$\boldsymbol{\mathcal{S}}\subseteq\boldsymbol{\mathcal{E}}$ such that each
$\boldsymbol{\mathcal{S}}_p:=\boldsymbol{\mathcal{S}}\cap\boldsymbol{\mathcal{E}}_p$ is a vector subspace and, for
each $p\in M$, there exist a trivialization
\[
\Phi\colon \boldsymbol{\mathcal{E}}\restriction_U\longrightarrow U\times F
\]
and a closed complemented subspace $F_0\subseteq F$ for which
\[
\Phi(\boldsymbol{\mathcal{S}}\restriction_U)=U\times F_0.
\]
\end{definition}

\begin{proposition}
\label{prop:subhaz-complementado-subvariedad-banach}
Let $M$ be a Banach manifold, and let $\boldsymbol{\mathcal E}\to M$ be a Banach vector bundle.
Every complemented subbundle $\boldsymbol{\mathcal{S}}\subseteq\boldsymbol{\mathcal{E}}$ is a Banach vector
bundle and, over each component on which its model fiber is fixed, its
inclusion in $\boldsymbol{\mathcal{E}}$ is a split embedded submanifold. In the
trivializations of Definition~\ref{def:subhaz-complementado-banach}, the local normal
bundle is modeled on any closed complement of $F_0$.
\end{proposition}

\begin{proof}
If two adapted trivializations intersect, their transition map
sends $F_0$ onto the corresponding subfiber model. The restriction is
a continuous linear isomorphism. To check regularity in norm, let $F_0$ and $\widetilde F_0$ be the subfiber models, and let
$\widetilde P\colon F\to\widetilde F_0$ be a continuous projection.
The restricted operator can be written as
$\widetilde P\circ g_{\beta\alpha}(p)\restriction_{F_0}$.
Restriction and composition with $\widetilde P$ are continuous linear
operations between the operator spaces; hence this family is $C^r$. Therefore, the restrictions of the trivializations define the
bundle structure of $\boldsymbol{\mathcal{S}}$. If $F=F_0\oplus F_1$, the inclusion is locally
\[
U\times F_0\longrightarrow U\times(F_0\oplus F_1),
\qquad
(p,v)\longmapsto(p,v,0).
\]
This is the local form of a split embedded submanifold.
\end{proof}

\begin{proposition}[Quotient by a complemented subbundle]
\label{prop:cociente-subhaz-complementado-banach}
Let $M$ be a Banach manifold, and let $\boldsymbol{\mathcal E}\to M$ be a Banach vector bundle.
If $\boldsymbol{\mathcal{S}}\subseteq\boldsymbol{\mathcal{E}}$ is a complemented subbundle, the union
\[
\boldsymbol{\mathcal{E}}/\boldsymbol{\mathcal{S}}
:=
\coprod_{p\in M}\boldsymbol{\mathcal{E}}_p/\boldsymbol{\mathcal{S}}_p
\]
restricted to each connected component $C\subseteq M$ has a natural
Banach vector bundle structure over $C$, with a fixed quotient model
$Q_C$. If the spaces $Q_C$ are all isomorphic to a single Banach
space $Q$, the full union is a bundle over $M$ with model fiber $Q$. On
each component, the sequence
\[
0\longrightarrow\boldsymbol{\mathcal{S}}
\longrightarrow\boldsymbol{\mathcal{E}}
\longrightarrow\boldsymbol{\mathcal{E}}/\boldsymbol{\mathcal{S}}
\longrightarrow0
\]
is exact fiberwise and splits locally by homomorphisms of class
$C^r$.
If $F=F_0\oplus F_1$ in an adapted trivialization, the projection
$F_1\longrightarrow F/F_0$ is a topological linear isomorphism. This identification
depends on the complement, whereas the quotient bundle is canonical.
\end{proposition}

\begin{proof}
A transition map between adapted trivializations sends the
subfiber model onto the corresponding model and therefore induces a
continuous linear isomorphism between the quotient spaces. These induced
operators depend on the base with $C^r$ regularity. Indeed, if
$q_\alpha\colon F\to F/F_{0,\alpha}$ is the quotient projection and
$j_\alpha\colon F/F_{0,\alpha}\to F$ a continuous linear section obtained from
a complement, the induced operator is
$q_\beta g_{\beta\alpha}(p)j_\alpha$. This expression is $C^r$ in the operator
norm. Its value does not depend on the chosen section, because two
representatives of a class differ by an element of $F_{0,\alpha}$, and the transition
map sends that subspace into $F_{0,\beta}$. The cocycle identities
follow by applying the operators to equivalence classes. On each connected component $C$, the Banach isomorphism type of
these quotient spaces is locally constant. Indeed, trivialization domains
that can be joined by a chain of nonempty intersections
have isomorphic quotient models, and the unions corresponding to the
different classes of chains would form a partition of $C$ into open sets;
connectedness forces there to be only one class. Moreover, $C$ is open because
a Banach manifold is locally connected. Fixing a model $Q_C$ and
topological linear identifications of the local models with $Q_C$, the
induced operators become a cocycle with values in
$\operatorname{GL}(Q_C)$.
Proposition~\ref{prop:transiciones-haz-banach} then constructs the quotient
bundle. If
$F=F_0\oplus F_1$, the canonical projection restricted to $F_1$ is bijective and
continuous; its inverse is continuous by the continuity of the projections
associated with the topological direct sum. This gives the local splitting.
\end{proof}

\begin{definition}[Section]
\label{def:seccion-haz-banach}
Let $M$ be a Banach manifold, and let $\boldsymbol{\mathcal E}\to M$ be a Banach vector bundle.
A \textbf{section} of $\boldsymbol{\mathcal E}\longrightarrow M$ over an open subset $U\subseteq M$ is
a map $\mathbf{s}\colon U\longrightarrow\boldsymbol{\mathcal E}$ such that $\pi\circ \mathbf{s}=\operatorname{id}_U$. It is of
class $C^r$ if, in every local trivialization, there is a map
$s_\alpha\in C^r(U\cap U_\alpha,F)$ for which
\[
\Phi_\alpha(\mathbf{s}(p))=(p,s_\alpha(p)).
\]
The vector space of global $C^r$ sections is denoted by
$\Gamma^r(\boldsymbol{\mathcal E})$.
\end{definition}

On an intersection we have
\[
s_\beta(p)=g_{\beta\alpha}(p)s_\alpha(p),
\]
so that the regularity of a section can be checked in any
family of trivializations.

\begin{proposition}[Local data and gluing of sections]
\label{prop:pegado-secciones-haz-banach}
Let $M$ be a Banach manifold, and let $\boldsymbol{\mathcal E}\to M$ be a Banach vector bundle of class $C^r$.
A family of maps $s_\alpha\in C^r(U_\alpha,F)$ determines a unique
section $\mathbf{s}\in\Gamma^r(\boldsymbol{\mathcal E})$ if and only if
\[
s_\beta(p)=g_{\beta\alpha}(p)s_\alpha(p)
\]
on each overlap. More generally, sections of class $C^r$ defined
over an open cover and agreeing pairwise on overlaps
glue to a unique section of class $C^r$. Moreover,
$\Gamma^r(\boldsymbol{\mathcal E})$ is a module over $C^r(M)$ under
$(f\mathbf{s})(p):=f(p)\mathbf{s}(p)$, and every $e\in\boldsymbol{\mathcal E}_p$ is the value at $p$ of a
local section of class $C^r$.
\end{proposition}

\begin{proof}
The compatibility condition allows us to define
\[
\mathbf{s}(p):=\Phi_\alpha^{-1}(p,s_\alpha(p))
\]
independently of the index $\alpha$. In each trivialization this section is
represented by $s_\alpha$, so it is of class $C^r$; the same formula
proves uniqueness. The argument for an arbitrary cover is identical.
The module operations are checked pointwise and preserve
regularity in the local representations. Finally, if
$\Phi\colon \boldsymbol{\mathcal E}\restriction_U\longrightarrow U\times F$ is a trivialization about $p$ and
$\Phi(e)=(p,v)$, the section $q\mapsto\Phi^{-1}(q,v)$ over $U$ takes the value
$e$ at $p$.
\end{proof}

\begin{proposition}[Pullback bundle]
\label{prop:pullback-haz-banach}
Let $M$ and $N$ be Banach manifolds.
Let $f\colon N\longrightarrow M$ be a map of class $C^r$, and let $\boldsymbol{\mathcal{E}}\longrightarrow M$ be a Banach
vector bundle of class $C^r$. The set
\[
f^*\boldsymbol{\mathcal{E}}
:=
\{(q,e)\in N\times\boldsymbol{\mathcal{E}}\mid f(q)=\pi(e)\},
\qquad
\pi_f(q,e):=q,
\]
has a natural Banach vector bundle structure over $N$, called
the \textbf{induced bundle} or \textbf{pullback}. Its transition maps are
\[
g_{\beta\alpha}\circ f
\]
over $f^{-1}(U_\alpha\cap U_\beta)$. This structure does not depend on the
chosen family of trivializations.
\end{proposition}

\begin{proof}
Over $f^{-1}(U_\alpha)$ define
\[
\widetilde\Phi_\alpha(q,e)
:=
\bigl(q,\operatorname{pr}_2\Phi_\alpha(e)\bigr).
\]
These bijections and their inverses provide the required local
trivializations. On an intersection, the fiber coordinate change is
$g_{\beta\alpha}(f(q))$, which is of class $C^r$. Compatibility with another
family of trivializations follows from
\eqref{eq:cociclos-equivalentes-haces-banach}: if the change is given by
$h_\alpha$, the induced cocycles are related by $h_\alpha\circ f$.
\end{proof}

\begin{proposition}[Direct sum]
\label{prop:suma-directa-haces-banach}
Let $M$ be a Banach manifold.
Let $\boldsymbol{\mathcal{E}}\longrightarrow M$ and $\boldsymbol{\mathcal{H}}\longrightarrow M$ be Banach vector bundles of
class $C^r$, with model fibers $F$ and $G$. The fiberwise union
\[
\boldsymbol{\mathcal{E}}\oplus\boldsymbol{\mathcal{H}}
:=
\coprod_{p\in M}(\boldsymbol{\mathcal{E}}_p\oplus\boldsymbol{\mathcal{H}}_p)
\]
is a Banach vector bundle with model fiber $F\oplus G$. In
simultaneous trivializations, its transition maps are
\[
g_{\beta\alpha}^{\boldsymbol{\mathcal{E}}}(p)
\oplus
g_{\delta\gamma}^{\boldsymbol{\mathcal{H}}}(p).
\]
The canonical inclusions and projections are bundle homomorphisms.
\end{proposition}

\begin{proof}
On $F\oplus G$ one may use, for example, the norm
$\|(v,w)\|=\|v\|_F+\|w\|_G$. The indicated diagonal operators are
invertible and depend on the base in a $C^r$ manner. The cocycle identities
hold componentwise.
\end{proof}

If $\boldsymbol{\mathcal E}$ and $\boldsymbol{\mathcal H}$ are bundles over the same base, we write
\[
\boldsymbol{\mathcal E}\times_M\boldsymbol{\mathcal H}
:=
\{(e,h)\in\boldsymbol{\mathcal E}\times\boldsymbol{\mathcal H}
\mid\pi_{\boldsymbol{\mathcal E}}(e)=\pi_{\boldsymbol{\mathcal H}}(h)\}
\]
for their fiber product. This notation describes the total space underlying
the direct sum when the linear operations are taken fiberwise.

\begin{proposition}[Dual bundle]
\label{prop:haz-dual-banach}\glsadd{haz-dual-topologico-banach}
Let $M$ and $N$ be Banach manifolds.
\index{continuous dual bundle}
Let $\boldsymbol{\mathcal{E}}\longrightarrow M$ be a Banach vector bundle of class $C^r$, with model
fiber $F$ and transition maps $g_{\beta\alpha}$. The union
\[
\boldsymbol{\mathcal{E}}'
:=
\coprod_{p\in M}\boldsymbol{\mathcal{E}}_p'
\]
has a natural Banach vector bundle structure of class $C^r$, with
model fiber $F'$. In the dual trivialization of $\Phi_\alpha$,
\[
\Phi_\alpha'(\lambda_p)
:=
\bigl(p,\lambda_p\circ\Phi_{\alpha,p}^{-1}\bigr),
\]
its transition maps are
\begin{equation}
\label{eq:transicion-haz-dual-banach}
g^\vee_{\beta\alpha}(p)
:=
\bigl(g_{\beta\alpha}(p)^{-1}\bigr)'.
\end{equation}
Fiberwise evaluation
\[
\operatorname{ev}\colon \boldsymbol{\mathcal{E}}'\times_M\boldsymbol{\mathcal{E}}\longrightarrow\mathbb R,
\qquad
(\lambda_p,e_p)\longmapsto\lambda_p(e_p),
\]
is of class $C^r$.

If $\mathscr A\colon \boldsymbol{\mathcal{E}}\longrightarrow\boldsymbol{\mathcal{H}}$ is a homomorphism over
$f\colon M\longrightarrow N$, there is a natural transpose homomorphism
\[
\mathscr A'\colon f^*\boldsymbol{\mathcal{H}}'\longrightarrow\boldsymbol{\mathcal{E}}',
\qquad
\mathscr A'_p(\mu_{f(p)}):=\mu_{f(p)}\circ\mathscr A_p.
\]
In particular, there is a canonical isomorphism
\[
\mathfrak d_f\colon (f^*\boldsymbol{\mathcal{H}})'\longrightarrow f^*(\boldsymbol{\mathcal{H}}').
\]
On the fiber over $p$, let
$\iota_p\colon \boldsymbol{\mathcal{H}}_{f(p)}\longrightarrow(f^*\boldsymbol{\mathcal{H}})_p$ be the
isomorphism $\iota_p(h):=(p,h)$. Then
\[
(\mathfrak d_f)_p(\lambda):=\lambda\circ\iota_p.
\]
Its inverse is $\mu\mapsto\mu\circ\iota_p^{-1}$. Both operators are linear and
continuous, and in the induced trivializations they are represented by the
identity on $G'$. Therefore, $\mathfrak d_f$ and its inverse are bundle
homomorphisms of class $C^r$.
\end{proposition}

\begin{proof}
Proposition~\ref{prop:dual-topologico-banach} shows that the operators
in \eqref{eq:transicion-haz-dual-banach} belong to
$\operatorname{GL}(F')$ and depend on $p$ with $C^r$ regularity. Moreover,
\[
g^\vee_{\gamma\beta}g^\vee_{\beta\alpha}
=
\bigl(g_{\gamma\beta}^{-1}\bigr)'
\bigl(g_{\beta\alpha}^{-1}\bigr)',
\]
and, since $(BA)'=A'B'$, the last expression is
\[
\bigl(g_{\beta\alpha}^{-1}g_{\gamma\beta}^{-1}\bigr)'
=
\bigl((g_{\gamma\beta}g_{\beta\alpha})^{-1}\bigr)'
=
g^\vee_{\gamma\alpha},
\]
so they satisfy the cocycle identity in the correct order.
Proposition~\ref{prop:transiciones-haz-banach} reconstructs the bundle. In the
trivializations of $\boldsymbol{\mathcal{E}}'$ and $\boldsymbol{\mathcal{E}}$, evaluation is the continuous
bilinear map $F'\times F\longrightarrow\mathbb R$; the dual transition law ensures that
these local expressions agree on overlaps.

The local representation of $\mathscr A'$ is
$A_{\beta\alpha}(p)'$. Transposition is continuous and linear in the operator
norm, and hence the family is of class $C^r$; its compatibility
follows by transposing
\eqref{eq:compatibilidad-homomorfismo-haces-banach}. Finally, the transition
maps of $(f^*\boldsymbol{\mathcal{H}})'$ and $f^*(\boldsymbol{\mathcal{H}}')$ are, in both cases,
\[
q\longmapsto
\bigl(h_{\delta\beta}(f(q))^{-1}\bigr)',
\]
so the fiberwise identity defines the canonical isomorphism.
\end{proof}

\begin{proposition}[Operator bundles]
\label{prop:haces-operadores-banach}
Let $M$ be a Banach manifold.
Let $\boldsymbol{\mathcal{E}}\longrightarrow M$ and $\boldsymbol{\mathcal{H}}\longrightarrow M$ be Banach vector bundles of
class $C^r$, with model fibers $F$ and $G$. The union
\[
\boldsymbol{\mathcal{L}}(\boldsymbol{\mathcal{E}},\boldsymbol{\mathcal{H}})
:=
\coprod_{p\in M}\mathcal{L}(\boldsymbol{\mathcal{E}}_p,\boldsymbol{\mathcal{H}}_p)
\]
is a Banach vector bundle with model fiber $\mathcal{L}(F,G)$. In two
pairs of trivializations, the coordinate change is
\begin{equation}
\label{eq:transicion-haz-operadores-banach}
A
\longmapsto
h_{\delta\beta}(p)A g_{\alpha\gamma}(p).
\end{equation}
Its sections of class $C^r$ are identified with the homomorphisms
$\boldsymbol{\mathcal{E}}\longrightarrow\boldsymbol{\mathcal{H}}$ of class $C^r$ over $\operatorname{id}_M$.

In particular,
\[
\operatorname{End}(\boldsymbol{\mathcal{E}})
:=
\boldsymbol{\mathcal{L}}(\boldsymbol{\mathcal{E}},\boldsymbol{\mathcal{E}})
\]
has fiber $\mathcal{L}(\boldsymbol{\mathcal{E}}_p)$ and transition maps
$A\mapsto g_{\beta\alpha}Ag_{\beta\alpha}^{-1}$. Evaluation
\[
\boldsymbol{\mathcal{L}}(\boldsymbol{\mathcal{E}},\boldsymbol{\mathcal{H}})\times_M\boldsymbol{\mathcal{E}}
\longrightarrow\boldsymbol{\mathcal{H}},
\qquad
(A_p,e_p)\longmapsto A_pe_p,
\]
is a bundle map of class $C^r$ and is bilinear in the fiber variables.
\end{proposition}

\begin{proof}
The space $\mathcal{L}(F,G)$ is Banach. Composition of operators is
continuous and bilinear, so the maps
\eqref{eq:transicion-haz-operadores-banach} are continuous linear
isomorphisms that depend on the point with $C^r$ regularity. The
cocycle identities follow from those of $g$ and $h$, and reconstruction follows from
Proposition~\ref{prop:transiciones-haz-banach}. The identification of
sections with homomorphisms is exactly the compatibility condition
\eqref{eq:compatibilidad-homomorfismo-haces-banach}. Evaluation is represented
locally by the continuous bilinear map
$\mathcal{L}(F,G)\times F\longrightarrow G$.
\end{proof}

\section{Tangent and cotangent bundles}
\label{sec:haces-tangente-cotangente-banach}

Let $M$ be a Banach manifold of class $C^k$, modeled on $X$, and
suppose first that $k\geq1$. Following
Definition~\ref{def:espacio-tangente-banach}, set
\[
TM:=\coprod_{p\in M}T_pM,
\qquad
\pi_{TM}(\xi_p):=p.
\]
Each chart $(U,\varphi)$ induces a trivialization
\[
\tau_\varphi\colon TM\restriction_U\longrightarrow U\times X,
\qquad
\tau_\varphi(\xi_p)
:=
\bigl(p,\Theta_{\varphi,p}(\xi_p)\bigr).
\]
If $(V,\psi)$ is another chart and $p\in U\cap V$, the change in the fiber is
\begin{equation}
\label{eq:transicion-haz-tangente-banach}
A_{\psi\varphi}(p)
:=
D(\psi\circ\varphi^{-1})(\varphi(p))
\in\operatorname{GL}(X).
\end{equation}
The chain rule gives the cocycle identities.

\begin{proposition}[Regularity of the tangent bundle]
\label{prop:regularidad-haz-tangente-banach}
If $M$ is a Banach manifold of class $C^k$, with
$2\leq k<\infty$, the preceding trivializations make $TM$ a Banach
vector bundle of class $C^{k-1}$. If $M$ is smooth, $TM$ is a smooth Banach
vector bundle. For a merely $C^1$ manifold, the same
construction produces a continuous vector bundle.
\end{proposition}

\begin{proof}
If the coordinate change $\psi\circ\varphi^{-1}$ is of class $C^k$, the
map assigning to each point its derivative, with values in
$\mathcal{L}(X)$ equipped with the operator norm, is of class $C^{k-1}$.
We consider the base with its underlying $C^{k-1}$ structure, described in
Remark~\ref{obs:atlas-regularidad-menor-banach}. The expression
\eqref{eq:transicion-haz-tangente-banach}, the openness of
$\operatorname{GL}(X)$, and reconstruction from cocycles prove the assertion. When $k=1$, this map is
continuous; when $k=\infty$, it is smooth.
\end{proof}

The loss of one derivative in the preceding proposition is part of the
structure: for $k<\infty$, one should not in general assign $C^k$ regularity to $TM$
on the basis of an atlas of class $C^k$.

\begin{definition}[Cotangent bundle]
\label{def:haz-cotangente-banach}\glsadd{cotangente-banach}
Let $M$ be a Banach manifold of class at least $C^1$.
\index{cotangent bundle!of a Banach manifold}
The \textbf{cotangent space} at $p$ is the continuous dual of the
kinematic tangent space,
\[
T'_pM:=(T_pM)'.
\]
The \textbf{cotangent bundle} is the union
\[
T'M:=\coprod_{p\in M}T'_pM,
\qquad
\pi_{T'M}(\alpha_p):=p.
\]
In the chart $(U,\varphi)$ it is trivialized by
\[
\tau_\varphi^*(\alpha_p)
=
\bigl(p,\alpha_p\circ\Theta_{\varphi,p}^{-1}\bigr)
\in U\times X'.
\]
\end{definition}

Thus, $T'M$ is the continuous dual bundle of $TM$ in the sense of
Proposition~\ref{prop:haz-dual-banach}. We use neither the algebraic dual of
$T_pM$ nor the dual of the operational tangent space. This choice ensures
that each fiber is a Banach space and that evaluation is continuous. The
prime distinguishes this bundle from the finite-dimensional cotangent bundle $T^*M$ and is
consistent with the notation $X'$ for the continuous dual of a Banach
space.

On an intersection, if $A_{\psi\varphi}$ is the operator in
\eqref{eq:transicion-haz-tangente-banach}, the cotangent transition map
is
\begin{equation}
\label{eq:transicion-haz-cotangente-banach}
\alpha
\longmapsto
\bigl(A_{\psi\varphi}(p)^{-1}\bigr)'\alpha.
\end{equation}
The map $A\mapsto A'$ is continuous and linear in the operator norm, and
inversion in $\operatorname{GL}(X)$ is smooth. Therefore, $T'M$ is a vector
bundle modeled on $X'$, of class $C^{k-1}$ when
$2\leq k<\infty$, continuous when $k=1$, and smooth when $M$ is smooth.

\begin{proposition}[Tangent--cotangent evaluation]
\label{prop:evaluacion-tangente-cotangente-banach}
Let $M$ be a Banach manifold of class at least $C^1$.
The bundle map
\[
\operatorname{ev}\colon T'M\times_M TM\longrightarrow\mathbb R,
\qquad
\operatorname{ev}(\alpha_p,\xi_p):=\alpha_p(\xi_p),
\]
is natural under coordinate changes. It has the same regularity
as the tangent and cotangent bundles; in particular, it is smooth if $M$ is smooth.
\end{proposition}

\begin{proof}
In a trivialization, the map is represented by continuous bilinear
evaluation
\[
X'\times X\longrightarrow\mathbb R,
\qquad
(\ell,v)\longmapsto\ell(v).
\]
Formula \eqref{eq:transicion-haz-cotangente-banach} was chosen precisely
to make this expression invariant.
\end{proof}

\begin{definition}[Covector fields and the differential]
\label{def:campo-covectorial-banach}
\index{covector field!on a Banach manifold}
\index{differential!of a function on a Banach manifold}
Let $M$ be a smooth Banach manifold. A \textbf{covector field} or a
\textbf{$1$-form} on an open subset $U\subseteq M$ is a smooth section
$\boldsymbol{\alpha}\colon U\longrightarrow T'M$. If $f\in C^\infty(U)$, its \textbf{differential} is the
covector field determined by
\[
(df)_p(v)
:=
\frac{d}{dt}\Biggr|_{t=0}f(\gamma(t)),
\qquad
v=[\gamma]\in T_pM.
\]
\end{definition}

The definition does not depend on the representative $\gamma$ because, in a chart
$(U,\varphi)$, we have
\begin{equation}
\label{eq:coordenadas-diferencial-banach}
(df)_\varphi(x)=D(f\circ\varphi^{-1})(x).
\end{equation}
In the same chart, a covector field is represented by
\[
\alpha_\varphi\colon \varphi(U)\longrightarrow X',
\qquad
\alpha_\varphi(\varphi(p))
:=
\boldsymbol{\alpha}_p\circ\Theta_{\varphi,p}^{-1},
\]
and, if $(V,\psi)$ is another chart, satisfies
\begin{equation}
\label{eq:transformacion-campo-covectorial-banach}
\alpha_\psi(\psi(p))
=
\bigl(A_{\psi\varphi}(p)^{-1}\bigr)'
\alpha_\varphi(\varphi(p)).
\end{equation}
Formula \eqref{eq:coordenadas-diferencial-banach} shows that $df$ is
smooth. More precisely, if $M$ and $f$ are of class $C^k$, with
$1\leq k<\infty$, then $T'M$ and the section $df$ are of class
$C^{k-1}$; the loss comes from differentiating the changes of charts and
$f$.

\begin{proposition}[Elementary properties of the differential]
\label{prop:propiedades-diferencial-funcion-banach}
Let $M$ be a smooth Banach manifold.
For $f,g\in C^\infty(M)$ we have
\[
d(fg)=f\,dg+g\,df.
\]
Moreover, every covector is locally the differential of a function at the
point: if $p\in M$ and $\lambda_p\in T'_pM$, there exist a neighborhood $U$ of $p$
and $f\in C^\infty(U)$ such that $(df)_p=\lambda_p$.
\end{proposition}

\begin{proof}
The Leibniz rule follows by applying the product rule to
\eqref{eq:coordenadas-diferencial-banach}. For the second assertion, let
$(U,\varphi)$ be a chart about $p$ and set
\[
\ell:=\lambda_p\circ\Theta_{\varphi,p}^{-1}\in X'.
\]
The function
$f(q):=\ell(\varphi(q)-\varphi(p))$ is smooth on $U$, and its local representation
has constant derivative $\ell$. Consequently, $(df)_p=\lambda_p$.
\end{proof}

\begin{definition}[Pullback of covector fields]
\label{def:pullback-covectores-banach}
\index{pullback!of a differential form on a Banach manifold}
Let $f\colon M\longrightarrow N$ be a smooth map between Banach manifolds. For a
$1$-form $\boldsymbol{\alpha}$ on $N$, define
\[
(f^*\boldsymbol{\alpha})_p
:=
\boldsymbol{\alpha}_{f(p)}\circ df_p\in T'_pM.
\]
The field $f^*\boldsymbol{\alpha}$ is called the \textbf{pullback} of $\boldsymbol{\alpha}$ by $f$.
\end{definition}

\begin{proposition}[Regularity and naturality of pullback]
\label{prop:pullback-covectores-banach}
Let $M$, $N$, and $P$ be smooth Banach manifolds, and let $f\colon M\to N$ be a smooth map.
The pullback of a smooth $1$-form is a smooth $1$-form. If
$g\colon N\longrightarrow P$ is smooth, then
\[
(g\circ f)^*=f^*\circ g^*,
\qquad
\operatorname{id}_M^*=\operatorname{id},
\qquad
d(h\circ f)=f^*(dh)
\]
for every $h\in C^\infty(N)$.
\end{proposition}

\begin{proof}
Let $\varphi$ and $\psi$ be charts in $M$ and $N$, respectively, and write
\[
\widehat f:=\psi\circ f\circ\varphi^{-1},
\qquad
\widehat\alpha:=\alpha_\psi.
\]
The local representation of the pullback is
\[
(f^*\boldsymbol{\alpha})_\varphi(x)
=
D\widehat f(x)'\widehat\alpha(\widehat f(x)).
\]
The map $x\mapsto D\widehat f(x)$ is smooth with respect to the operator
norm, transposition is continuous and linear, and evaluation is continuous
and bilinear. Therefore, the expression is smooth. The three identities follow
from the chain rule and the pointwise definition.
\end{proof}

To treat forms of higher degree without choosing a tensor completion,
if $m\in\mathbb N$ and $X_0,Y_0$ are Banach spaces, denote by
\[
\mathcal{L}^m_{\mathrm{alt}}(X_0;Y_0)
\]
the subspace of continuous $m$-linear maps
$A\colon X_0^m\longrightarrow Y_0$ that are alternating, equipped with
the operator norm. It is a closed subspace of
$\mathcal{L}^m(X_0;Y_0)$ and therefore is a Banach space.
Coordinate changes act by
\[
\omega
\longmapsto
\omega\circ(A^{-1},\dots,A^{-1}).
\]
Proposition~\ref{prop:transiciones-haz-banach} thus provides a Banach
bundle, denoted by $\bigwedge^m(T'M)$, whose fiber at $p$ is
$\mathcal{L}^m_{\mathrm{alt}}(T_pM;\mathbb R)$.

\begin{definition}[Differential forms]
\label{def:formas-diferenciales-banach}
Let $M$ and $N$ be smooth Banach manifolds.
\index{differential form!on a Banach manifold}
A \textbf{smooth differential $m$-form} on $M$ is a smooth section of
$\bigwedge^m(T'M)$, for $m\in\mathbb N$. The space of these sections is
denoted by $\Omega^m(M)$; for $m=0$ we set $\Omega^0(M):=C^\infty(M)$.
If $f\colon M\longrightarrow N$ is smooth and $\boldsymbol{\omega}\in\Omega^m(N)$, its pullback is defined by
\begin{equation}
\label{eq:pullback-formas-banach}
(f^*\boldsymbol{\omega})_p(v_1,\dots,v_m)
:=
\boldsymbol{\omega}_{f(p)}(df_pv_1,\dots,df_pv_m).
\end{equation}
\end{definition}

The same local proof as in
Proposition~\ref{prop:pullback-covectores-banach}, using continuous multilinear
composition, shows that \eqref{eq:pullback-formas-banach} defines a
smooth $m$-form and that pullback is functorial. This formulation does not identify
continuous forms with a completed tensor product; it thus avoids imposing
a tensor norm that is not canonical.

Definitions~\ref{def:haz-cotangente-banach} and
\ref{def:campo-covectorial-banach}, together with Proposition
\ref{prop:propiedades-diferencial-funcion-banach}, correspond to
Definitions~\ref{def:nociones-fundamentales-haz-cotangente-de} and
\ref{def:nociones-fundamentales-campo-covectorial} and
Corollary~\ref{cor:nociones-fundamentales-diferencial-funcion-suave-campo-covectorial-suave}.
The pullback formula is the same as in
Definition~\ref{def:nociones-fundamentales-pullback-de-bajo-en} and preserves
the identities of Proposition~\ref{propiedades del pullback}. The functional-analytic difference
is that here each covector belongs to the continuous dual and all
regularity is checked in the operator norm; no finite-dimensional
identification with a tensor product is used.

\begin{remark}[Bidual and tensor products]
\label{obs:doble-dual-tensores-banach}
There is a canonical isometric linear injection
\[
J_X\colon X\longrightarrow X'',
\qquad
J_X(v)(\ell):=\ell(v),
\]
but it need not be surjective. Hence, if the model space is not reflexive,
there is no canonical fiberwise identification between $TM$ and $(T'M)'$.

Nor is there a unique tensor completion implicit in the Banach
category. For example, $X'\otimes F$ maps to the space of finite-rank operators
from $X$ to $F$, whereas its projective and injective completions
may differ and do not in general coincide with $\mathcal L(X,F)$.
Consequently, tensor bundles must be defined by specifying a complete
tensor norm. When only continuous multilinear maps
are required, it is canonical to work directly with fibers such as
\[
\mathcal L^m(T_pM;\mathbb R)
\quad\text{or}\quad
\mathcal L(T_pM,T_pM),
\]
without imposing a tensor identification that may fail.
\end{remark}

\section{Vector fields and integral curves}
\label{sec:campos-curvas-integrales-banach}

Throughout the rest of this chapter, $M$ will be a smooth Banach manifold modeled
on a Banach space $X$.

\begin{definition}[Vector field]
\label{def:campo-vectorial-banach}
Let $M$ be a smooth Banach manifold.
\index{vector field!on a Banach manifold}
A \textbf{vector field} on $M$ is a section
$\mathbf{X}\colon M\longrightarrow TM$. In a chart $(U,\varphi)$ it is represented by
\[
X_\varphi\colon \varphi(U)\longrightarrow X,
\qquad
X_\varphi(\varphi(p)):=\Theta_{\varphi,p}(\mathbf{X}(p)).
\]
On the overlap of two charts we have
\begin{equation}
\label{eq:transformacion-campo-vectorial-banach}
X_\psi(\psi(p))
=
D(\psi\circ\varphi^{-1})(\varphi(p))X_\varphi(\varphi(p)).
\end{equation}
The field is continuous, locally Lipschitz, or of class $C^r$ if its
coordinate representations have the corresponding property.
\end{definition}

We denote by $\mathfrak X(U):=\Gamma^\infty(TU)$ the space of smooth
vector fields on an open subset $U\subseteq M$. If
$f\in C^\infty(U)$ and $\mathbf{X}\in\mathfrak X(U)$, the function
$\mathbf{X}(f)\in C^\infty(U)$ is defined by
\[
\mathbf{X}(f)(p):=(df)_p(\mathbf{X}(p)),
\qquad p\in U.
\]
If
$\mathbf{X},\mathbf{Y}\in\mathfrak X(U)$, in a chart we define
\begin{equation}
\label{eq:corchete-campos-banach}
[\mathbf{X},\mathbf{Y}]_\varphi(x)
:=
DY_\varphi(x)X_\varphi(x)-DX_\varphi(x)Y_\varphi(x).
\end{equation}

\begin{proposition}[Bracket of vector fields]
\label{prop:corchete-campos-banach}
Let $M$ be a smooth Banach manifold, and let $U\subseteq M$ be open.
Expression \eqref{eq:corchete-campos-banach} is independent of the chart and
defines a smooth vector field $[\mathbf{X},\mathbf{Y}]$. Moreover,
\[
[\mathbf{X},\mathbf{Y}]=-[\mathbf{Y},\mathbf{X}],
\qquad
[\mathbf{X},f\mathbf{Y}]=\mathbf{X}(f)\mathbf{Y}+f[\mathbf{X},\mathbf{Y}]
\]
for every $f\in C^\infty(U)$, and the bracket satisfies the Jacobi identity.
\end{proposition}

\begin{proof}
Let $\chi$ be a coordinate change, and write
$\widetilde X(\chi(x))=D\chi(x)X(x)$ and
$\widetilde Y(\chi(x))=D\chi(x)Y(x)$. Differentiating and subtracting gives
\begin{align*}
D\widetilde Y(\chi(x))\widetilde X(\chi(x))
&-D\widetilde X(\chi(x))\widetilde Y(\chi(x))\\
={}&D\chi(x)\bigl(DY(x)X(x)-DX(x)Y(x)\bigr),
\end{align*}
because the terms containing $D^2\chi(x)(X(x),Y(x))$ cancel by the
symmetry of the second derivative. This is exactly the transformation
law for a vector field. The first two identities are
checked in coordinates using the product rule. For Jacobi,
first consider the action on a function $f\in C^\infty(U)$.
If $\widehat f=f\circ\varphi^{-1}$, the product rule gives
\[
 \mathbf X(\mathbf Y(f))
 =D^2\widehat f(X_\varphi,Y_\varphi)
  +D\widehat f\bigl(DY_\varphi[X_\varphi]\bigr)
\]
in the chart coordinates. The symmetry of $D^2\widehat f$ implies
$[\mathbf X,\mathbf Y](f)=\mathbf X(\mathbf Y(f))-\mathbf Y(\mathbf X(f))$.
Applying this identity twice, the three terms of Jacobi act as
\begin{align*}
 [\mathbf X,[\mathbf Y,\mathbf Z]](f)
 &=\mathbf X\mathbf Y\mathbf Z(f)-\mathbf X\mathbf Z\mathbf Y(f)
   -\mathbf Y\mathbf Z\mathbf X(f)+\mathbf Z\mathbf Y\mathbf X(f),\\
 [\mathbf Y,[\mathbf Z,\mathbf X]](f)
 &=\mathbf Y\mathbf Z\mathbf X(f)-\mathbf Y\mathbf X\mathbf Z(f)
   -\mathbf Z\mathbf X\mathbf Y(f)+\mathbf X\mathbf Z\mathbf Y(f),\\
 [\mathbf Z,[\mathbf X,\mathbf Y]](f)
 &=\mathbf Z\mathbf X\mathbf Y(f)-\mathbf Z\mathbf Y\mathbf X(f)
   -\mathbf X\mathbf Y\mathbf Z(f)+\mathbf Y\mathbf X\mathbf Z(f).
\end{align*}
Each composition appears once with each sign, so the sum is zero.
At each point $p$, every continuous covector is the differential of a local
function by Proposition~\ref{prop:propiedades-diferencial-funcion-banach}.
Since the continuous dual separates points by Hahn--Banach, a vector annihilated by
all these differentials is zero. The sum of the three fields is therefore
the zero field, which is the Jacobi identity.
\end{proof}

The local Lipschitz condition is independent of the chart. Indeed, let
$\chi$ be a coordinate change, and let $B$ be a convex ball on which
$X_\varphi$ is $L_X$-Lipschitz. After shrinking $B$, the smoothness of
$\chi$ provides constants $C_0,C_1,C_2>0$ such that, for $x,y\in B$,
\[
\|D\chi(x)\|\leq C_0,
\qquad
\|D\chi(x)-D\chi(y)\|\leq C_1\|x-y\|,
\qquad
\|x-y\|\leq C_2\|\chi(x)-\chi(y)\|.
\]
If also $\|X_\varphi(x)\|\leq C_X$ on $B$, the transformation law gives
\[
\|X_\psi(\chi(x))-X_\psi(\chi(y))\|
\leq
C_2(C_0L_X+C_1C_X)\|\chi(x)-\chi(y)\|.
\]
This transports a Lipschitz constant between charts. In particular,
every field of class $C^1$ is locally Lipschitz, since its derivative is
bounded on a sufficiently small ball.

\begin{definition}[Integral curve]
\label{def:curva-integral-banach}
Let $M$ be a smooth Banach manifold, and let $\mathbf{X}$ be a vector field on $M$.
\index{integral curve!on a Banach manifold}
If $I\subseteq\mathbb R$ is an open interval and
$\gamma\colon I\longrightarrow M$ is a curve of class $C^1$, its velocity at $t$ is
\[
\dot\gamma(t):=d\gamma_t(1)\in T_{\gamma(t)}M.
\]
The curve $\gamma$ is an \textbf{integral curve} of $\mathbf{X}$ if
\[
\dot\gamma(t)=\mathbf{X}(\gamma(t))
\]
for every $t\in I$. It is \textbf{maximal} if it is not the restriction of an integral
curve defined on a strictly larger interval.
\end{definition}

In a chart, the integral curve equation becomes the autonomous
ODE
\[
x'(t)=X_\varphi(x(t)).
\]
The transformation law \eqref{eq:transformacion-campo-vectorial-banach} and the
chain rule ensure that this equation is independent of the coordinates.

\begin{theorem}[Existence and uniqueness of integral curves]
\label{teo:existencia-curvas-integrales-banach}
Let $M$ be a smooth Banach manifold.
Let $\mathbf{X}$ be a continuous, locally Lipschitz field on $M$. For each
$p\in M$ there exists a unique maximal integral curve
\[
\gamma_p\colon I_p\longrightarrow M,
\qquad
0\in I_p,
\qquad
\gamma_p(0)=p.
\]
The domain $I_p$ is an open interval. Every integral curve taking the
value $p$ at time zero is a restriction of $\gamma_p$.
\end{theorem}

\begin{proof}
Theorem~\ref{teo:picard-lindelof-banach}, applied to a coordinate
representation of $\mathbf{X}$, gives a local integral curve.
Proposition~\ref{prop:unicidad-soluciones-edo-banach} shows that two integral
curves with a common value agree near that time and, by the
open-and-closed argument on the intersection of their intervals,
agree on every common component. The local solutions can therefore
be glued. The same union argument as in
Theorem~\ref{teo:solucion-maximal-edo-banach} produces a maximal open
interval and a unique integral curve on it.
\end{proof}

\section{Flows on Banach manifolds}
\label{sec:flujos-variedades-banach}

\begin{definition}[Local flow]
\label{def:flujo-local-banach}
Let $M$ be a smooth Banach manifold.
\index{local flow!on a Banach manifold}
A \textbf{local flow} on $M$ is a continuous map
\[
\Phi\colon \mathcal D\longrightarrow M,
\]
where $\mathcal D\subseteq\mathbb R\times M$ is open and contains
$\{0\}\times M$. For each $p\in M$, we require that
\[
 \mathcal D_p:=\{t\in\mathbb R\mid(t,p)\in\mathcal D\}
\]
be an interval containing $0$. The map satisfies
\[
\Phi(0,p)=p
\]
and, if $s\in\mathcal D_p$, the exact equivalence of domains
\[
 t\in\mathcal D_{\Phi(s,p)}
 \quad\Longleftrightarrow\quad
 t+s\in\mathcal D_p,
\]
together with the law
\[
\Phi(t+s,p)=\Phi(t,\Phi(s,p))
\]
whenever those conditions hold. If
$\mathcal D=\mathbb R\times M$, it is called a \textbf{global flow}. A flow
is of class $C^r$, with $r\in\mathbb N_0\cup\{\infty\}$, if the map $\Phi$ is
of that class. It is said to have
\textbf{generating field} $\mathbf{X}$ if each curve $t\mapsto\Phi(t,p)$ is of class
$C^1$ and
\[
\frac{\partial}{\partial t}\Biggr|_{t=0}\Phi(t,p)=\mathbf{X}(p)
\]
for every $p\in M$. A local flow with generating field $\mathbf{X}$ is
\textbf{maximal} if every other local flow with generator $\mathbf{X}$ is a restriction
of it. In particular, it admits no proper extension with the same generator.
\end{definition}

For a field $\mathbf{X}$ as in
Theorem~\ref{teo:existencia-curvas-integrales-banach}, define
\[
\mathcal D_{\mathbf{X}}
:=
\{(t,p)\in\mathbb R\times M\mid t\in I_p\},
\qquad
\Phi_{\mathbf{X}}(t,p):=\gamma_p(t).
\]

\begin{theorem}[Maximal flow and smooth dependence]
\label{teo:flujo-maximal-banach}
Let $M$ be a smooth Banach manifold.
Let $\mathbf{X}$ be a continuous, locally Lipschitz field on $M$. Then
$\mathcal D_{\mathbf{X}}$ is open, $\Phi_{\mathbf{X}}$ is continuous, and it is the unique maximal
local flow whose orbit curves are integral curves of $\mathbf{X}$.

If $\mathbf{X}$ is of class $C^r$, with
$r\in\mathbb N\cup\{\infty\}$, then $\Phi_{\mathbf{X}}$ is of class $C^r$. For
any $p\in M$ and $s,t\in\mathbb R$ such that
$(s,p)\in\mathcal D_{\mathbf{X}}$, we have the equivalence
\[
(t,\Phi_{\mathbf{X}}(s,p))\in\mathcal D_{\mathbf{X}}
\quad\Longleftrightarrow\quad
(t+s,p)\in\mathcal D_{\mathbf{X}},
\]
and, when these conditions hold,
\begin{equation}
\label{eq:ley-grupo-flujo-banach}
\Phi_{\mathbf{X}}(t,\Phi_{\mathbf{X}}(s,p))=\Phi_{\mathbf{X}}(t+s,p).
\end{equation}
\end{theorem}

\begin{proof}
Theorem~\ref{teo:dependencia-parametros-edo-banach}, with the initial
point as parameter, gives a local solution map around each $q\in M$,
\[
 \Psi\colon(-a,a)\times N\longrightarrow M,
 \qquad a>0,\quad q\in N.
\]
It is continuous, and of class $C^r$ when $\mathbf X$ is of class
$C^r$. Its curves are integral curves of $\mathbf X$, and
$\Psi(0,p)=p$.

Let $(t_0,p_0)\in\mathcal D_{\mathbf X}$ with $t_0>0$.
The image of $\gamma_{p_0}|_{[0,t_0]}$ is compact. Cover it by
finitely many domains $N$ of this kind and choose $a_*>0$ smaller
than all their existence times. Choose a partition
\[
 0=s_0<s_1<\cdots<s_m=t_0,
 \qquad \Delta_j:=s_{j+1}-s_j<a_*\quad(0\leq j<m).
\]
For each $j$, choose a domain $N_j$ from the cover containing
$\gamma_{p_0}(s_j)$, and denote its local map by
$\Psi_j\colon I_j\times N_j\to M$, where $I_j$ is an open
interval containing $[0,\Delta_j]$.
Define recursively
\[
 V_0=M,\quad F_0=\operatorname{id}_M,\qquad
 V_{j+1}=V_j\cap F_j^{-1}(N_j),\quad
 F_{j+1}(p)=\Psi_j(\Delta_j,F_j(p)).
\]
The preimage is taken in the domain $V_j$ of $F_j$.
Inductively, $V_j$ is open, contains $p_0$, and $F_j$ has the
asserted regularity. Uniqueness of integral curves gives
$F_j(p_0)=\gamma_{p_0}(s_j)$, which justifies
$p_0\in V_{j+1}$ at each step.

Choose $0<\varepsilon<\Delta_{m-1}$ such that
$(\Delta_{m-1}-\varepsilon,\Delta_{m-1}+\varepsilon)\subset I_{m-1}$.
For $(t,p)\in(t_0-\varepsilon,t_0+\varepsilon)\times V_m$,
concatenating the first $m-1$ segments with the variable last segment
gives an integral curve, and uniqueness yields
\[
 \Phi_{\mathbf X}(t,p)
 =\Psi_{m-1}\bigl(t-s_{m-1},F_{m-1}(p)\bigr).
\]
Thus this open product lies in $\mathcal D_{\mathbf X}$, and the
formula proves continuity and regularity there. For $t_0<0$, apply
the construction to $-\mathbf X$; at $t_0=0$, use the local map.
This proves the assertions throughout the maximal domain.

Now fix $s$ with $(s,p)\in\mathcal D_{\mathbf{X}}$. The curves
\[
t\longmapsto\Phi_{\mathbf{X}}(t,\Phi_{\mathbf{X}}(s,p))
\quad\text{and}\quad
t\longmapsto\Phi_{\mathbf{X}}(t+s,p)
\]
are integral curves of $\mathbf{X}$ and take the same value at $t=0$. Uniqueness
shows that they agree on their common domain. Maximality of both curves
identifies their intervals of definition and proves both the equivalence of
domains and \eqref{eq:ley-grupo-flujo-banach}. Finally, any flow
with generating field $\mathbf{X}$ consists of integral curves: differentiating
$\Phi(t+h,p)=\Phi(h,\Phi(t,p))$ at $h=0$ gives
$\partial_t\Phi(t,p)=\mathbf X(\Phi(t,p))$. Uniqueness of integral
curves then shows that this flow is a restriction of
$\Phi_{\mathbf{X}}$.
\end{proof}

The preceding notions extend
Definitions~\ref{def:nociones-fundamentales-campo-vectorial-en},
\ref{def:nociones-fundamentales-integral-variedad-suave-frontera-curva}, and
\ref{def:nociones-fundamentales-flujo}. The proof of
Theorem~\ref{teo: fundamental de flujos} retains its uniqueness and
gluing argument here; the new step is that local existence and dependence
on the initial value come from ODEs in Banach spaces,
Theorems~\ref{teo:picard-lindelof-banach} and
\ref{teo:dependencia-parametros-edo-banach}, without using local compactness.

For each $t\in\mathbb R$, let
\[
M_t:=\{p\in M\mid(t,p)\in\mathcal D_{\mathbf{X}}\}.
\]
The openness of $\mathcal D_{\mathbf{X}}$ implies that $M_t$ is open. The flow law
shows that
\[
\Phi_{\mathbf{X}}^t:=\Phi_{\mathbf{X}}(t,\cdot)\colon M_t\longrightarrow M_{-t}
\]
is a homeomorphism with inverse $\Phi_{\mathbf{X}}^{-t}$. If $\mathbf{X}$ is of class $C^r$,
both maps are of class $C^r$, and hence $\Phi_{\mathbf{X}}^t$ is a $C^r$ diffeomorphism
between open sets.

\begin{remark}[Fields and flows with $C^{1,1}$ regularity]
\label{obs:campos-flujos-atlas-C11}
Let $M$ be a Banach manifold with the $C^{1,1}$ atlas of
Definition~\ref{def:variedad-banach-C11}. A vector field is a section
$\mathbf{X}\colon M\longrightarrow TM$, and we say it is locally
Lipschitz if its representations in the charts of that atlas are locally Lipschitz.
This property does not depend on the chosen chart. Indeed, the law
\eqref{eq:transformacion-campo-vectorial-banach} uses only first
derivatives. If $\chi=\psi\circ\varphi^{-1}$, we can shrink the domain
about the point under consideration so that $D\chi$ is bounded and
Lipschitz, $\mathbf{X}_\varphi$ is bounded and Lipschitz, and $\chi^{-1}$
is Lipschitz on a ball in its domain. The last property follows
from the continuity of $D\chi^{-1}$ and the mean value inequality.
The estimate in Section~\ref{sec:campos-curvas-integrales-banach}
applies with these same bounds and shows that
$\mathbf{X}_\psi=(D\chi\cdot\mathbf{X}_\varphi)\circ\chi^{-1}$ is
locally Lipschitz.

The definition of an integral curve and the conclusions concerning existence,
uniqueness, and a continuous maximal flow in
Theorems~\ref{teo:existencia-curvas-integrales-banach} and
\ref{teo:flujo-maximal-banach} remain valid at this regularity.
To see this, in each chart we apply
Theorem~\ref{teo:picard-lindelof-banach} to the coordinate field. The chain
rule for $C^1$ changes transforms a coordinate solution
into an integral curve independent of the chart. Local uniqueness
allows these solutions to be glued and the union of their intervals to be taken, as
in Theorem~\ref{teo:existencia-curvas-integrales-banach}. By
Theorem~\ref{teo:dependencia-parametros-edo-banach}, the local solution
maps are continuous in time and in the initial value.
Composition over a finite partition of a compact segment of the
orbit, constructed in the proof of
Theorem~\ref{teo:flujo-maximal-banach}, proves openness of the flow
domain and joint continuity. Uniqueness gives the group
law and the time homeomorphisms together with their inverses.
No differentiability of the flow with respect to the initial value is asserted here
for a field that is only locally Lipschitz.
\end{remark}

\begin{proposition}[Generating field]
\label{prop:campo-generador-flujo-banach}
Let $M$ be a smooth Banach manifold.
If $\Phi\colon \mathcal D\longrightarrow M$ is a local flow of class $C^1$, then
\[
\mathbf{X}_\Phi(p)
:=
\frac{\partial}{\partial t}\Biggr|_{t=0}\Phi(t,p)
\in T_pM
\]
defines a continuous vector field. Each curve
$t\mapsto\Phi(t,p)$ is an integral curve of $\mathbf{X}_\Phi$. If $\Phi$ is of class
$C^r$, with $r\geq1$, then $\mathbf{X}_\Phi$ is of class $C^{r-1}$; if $\Phi$ is
smooth, so is $\mathbf{X}_\Phi$.
\end{proposition}

\begin{proof}
The identity $\Phi(0,p)=p$ shows that $\mathbf{X}_\Phi(p)\in T_pM$. Differentiating the
flow law
\[
\Phi(t+h,p)=\Phi(h,\Phi(t,p))
\]
with respect to $h$ at $h=0$ gives
\[
\frac{\partial}{\partial t}\Phi(t,p)
=
\mathbf{X}_\Phi(\Phi(t,p)).
\]
The regularity assertion follows by taking a partial derivative with respect
to time.
\end{proof}

\begin{theorem}[Global flow of a bounded field on a Banach space]
\label{teo:flujo-global-campo-acotado-banach}
Let $X$ be a Banach space, regarded as a smooth Banach manifold, and let $\mathbf{V}\colon X\longrightarrow X$ be a continuous, locally
Lipschitz field. If there exists $C\geq0$ such that
\[
\|\mathbf{V}(x)\|_X\leq C
\]
for every $x\in X$, then its maximal flow is global:
\[
\Phi_{\mathbf{V}}\colon \mathbb R\times X\longrightarrow X.
\]
Moreover,
\[
\|\Phi_{\mathbf{V}}(t,x)-x\|_X\leq C|t|
\]
for every $(t,x)\in\mathbb R\times X$. If $\mathbf{V}$ is of class $C^r$, with
$r\in\mathbb N\cup\{\infty\}$, then $\Phi_{\mathbf{V}}$ is of class $C^r$, and each
$\Phi_{\mathbf{V}}^t\colon X\longrightarrow X$ is a diffeomorphism of class $C^r$, with inverse
$\Phi_{\mathbf{V}}^{-t}$.
\end{theorem}

\begin{proof}
Corollary~\ref{cor:existencia-global-campo-acotado-edo-banach}, applied to
$F(t,x):=\mathbf{V}(x)$, shows that every maximal integral curve is defined for
all $t\in\mathbb R$ and provides the indicated estimate. Regularity of the
flow follows from
Theorem~\ref{teo:dependencia-parametros-edo-banach}. The law
\eqref{eq:ley-grupo-flujo-banach} gives
\[
\Phi_{\mathbf{V}}^{-t}\circ\Phi_{\mathbf{V}}^t
=
\operatorname{id}_X
=
\Phi_{\mathbf{V}}^t\circ\Phi_{\mathbf{V}}^{-t}.
\]
\end{proof}

The boundedness condition in the preceding theorem uses the distinguished norm of the space
$X$. On an abstract Banach manifold it is not an invariant condition until
additional structure is fixed, for example a continuous norm on each
fiber of $TM$ with suitable uniform control. The basic intrinsic criterion for
continuation remains the analysis, in a final chart, of convergence of the
curve and membership of its limit in the domain, as in
Theorem~\ref{teo:criterio-prolongacion-edo-banach}.

\section{Linear connections on Banach bundles}
\label{sec:conexiones-haces-banach}

Let $\pi\colon \boldsymbol{\mathcal{E}}\longrightarrow M$ be a smooth Banach vector bundle with model fiber
$F$, over a smooth Banach manifold modeled on $X$. In infinite
dimensions it is not enough to require algebraic identities of an operator on
sections: one must also require its local coefficients to depend
smoothly on the point in the operator norm topology.

\begin{definition}[Linear connection]
\label{def:conexion-lineal-haz-banach}
Let $M$ be a smooth Banach manifold, and let $\boldsymbol{\mathcal E}\to M$ be a smooth Banach vector bundle.
\index{linear connection!on a Banach bundle}
A \textbf{linear connection} on $\boldsymbol{\mathcal{E}}$ is a rule that, for each
open set $U\subseteq M$, each $\mathbf{X}\in\mathfrak X(U)$, and each
$\mathbf{s}\in\Gamma^\infty(\boldsymbol{\mathcal{E}}\restriction_U)$, assigns a section
\[
\nabla_{\mathbf{X}}\mathbf{s}\in\Gamma^\infty(\boldsymbol{\mathcal{E}}\restriction_U),
\]
compatibly with restriction to smaller open sets, and satisfies
\begin{align*}
\nabla_{f\mathbf{X}+g\mathbf{Y}}\mathbf{s}
&=f\nabla_{\mathbf{X}}\mathbf{s}+g\nabla_{\mathbf{Y}}\mathbf{s},\\
\nabla_{\mathbf{X}}(a\mathbf{s}+b\mathbf{t})
&=a\nabla_{\mathbf{X}}\mathbf{s}+b\nabla_{\mathbf{X}}\mathbf{t},\\
\nabla_{\mathbf{X}}(f\mathbf{s})
&=\mathbf{X}(f)\mathbf{s}+f\nabla_{\mathbf{X}}\mathbf{s}
\end{align*}
for $f,g\in C^\infty(U)$ and $a,b\in\mathbb R$. The following local regularity
is also required. If $(U,\varphi)$ is a chart of $M$ and
$\Phi\colon \boldsymbol{\mathcal{E}}\restriction_U\longrightarrow U\times F$ is a trivialization, there is a map
smooth in the operator norm,
\[
\Gamma_{\varphi,\Phi}\colon \varphi(U)\longrightarrow
\mathcal{L}\bigl(X,\mathcal{L}(F)\bigr)
\]
such that
\begin{equation}
\label{eq:conexion-local-haz-banach}
\bigl(\nabla_{\mathbf{X}}\mathbf{s}\bigr)_{\varphi,\Phi}(x)
=
Ds_{\varphi,\Phi}(x)[X_\varphi(x)]
+
\Gamma_{\varphi,\Phi}(x)(X_\varphi(x))s_{\varphi,\Phi}(x).
\end{equation}
The map $\Gamma_{\varphi,\Phi}$ is called the \textbf{local coefficient of
the connection} in the chosen chart and trivialization.
\end{definition}

The first identity expresses linearity over $C^\infty(U)$ in the
direction, whereas the third is the Leibniz rule in the section. The
additional local condition ensures that the operation is a continuous
differential structure and not merely an algebraic operator.

The algebraic part agrees with Definition
\ref{definicion 1 conexion}. The difference is analytic: in a
trivialization with infinite-dimensional fiber, the coefficient must be a smooth map with
values in $\mathcal L(X,\mathcal L(F))$ equipped with the operator norm.
This hypothesis replaces the automatic control of matrix families in
finite dimensions and ensures smoothness of the evaluation appearing in
\eqref{eq:conexion-local-haz-banach}.

\begin{proposition}[Transformation law for a connection]
\label{prop:transformacion-conexion-haz-banach}
Let $M$ be a smooth Banach manifold, and let $\boldsymbol{\mathcal E}\to M$ be a smooth Banach vector bundle with connection $\nabla$.
Let $(U,\varphi,\Phi)$ and $(V,\psi,\Psi)$ be two choices of chart and
trivialization. On the overlap, write
\[
\chi:=\psi\circ\varphi^{-1},
\qquad
g(x):=g_{\Psi\Phi}(\varphi^{-1}(x)).
\]
If $\Gamma$ and $\widetilde\Gamma$ are the corresponding local
coefficients, then
\begin{equation}
\label{eq:ley-transformacion-conexion-banach}
\widetilde\Gamma(\chi(x))(D\chi(x)v)
=
g(x)\Gamma(x)(v)g(x)^{-1}
-Dg(x)[v]g(x)^{-1}
\end{equation}
for every $v\in X$. Conversely, a family of smooth maps with values
in $\mathcal L(X,\mathcal L(F))$ satisfying this law determines a unique
linear connection.
\end{proposition}

\begin{proof}
Let $x=\varphi(p)$, $w=s_{\varphi,\Phi}(x)$, and
$v=X_\varphi(x)$. In the second coordinates we have
\[
s_{\psi,\Psi}(\chi(x))=g(x)w,
\qquad
X_\psi(\chi(x))=D\chi(x)v.
\]
The chain rule and product rule give
\[
D(s_{\psi,\Psi}\circ\chi)(x)[v]
=
Dg(x)[v]w+g(x)Ds_{\varphi,\Phi}(x)[v].
\]
Since $\nabla_{\mathbf{X}}\mathbf{s}$ is a section of $\boldsymbol{\mathcal E}$, its second representation
must be $g(x)$ times the first. Substituting
\eqref{eq:conexion-local-haz-banach} and canceling the term involving
$Ds_{\varphi,\Phi}$ gives
\[
\widetilde\Gamma(\chi(x))(D\chi(x)v)g(x)w
=
g(x)\Gamma(x)(v)w-Dg(x)[v]w.
\]
Since $w\in F$ is arbitrary, this identity is equivalent to
\eqref{eq:ley-transformacion-conexion-banach}.

Conversely, define $\nabla_{\mathbf{X}}\mathbf{s}$ by
\eqref{eq:conexion-local-haz-banach} in each trivialization. The preceding
calculation, read in reverse, shows that the local expressions
transform as a section and therefore glue together. The identities in
Definition~\ref{def:conexion-lineal-haz-banach} are checked directly, and
the local formula proves uniqueness.
\end{proof}

\begin{theorem}[Existence of connections]
\label{teo:existencia-conexiones-haces-banach}
Let $M$ be a smooth Banach manifold.
Suppose that $M$ admits smooth partitions of unity subordinate to every
open cover. Then every smooth Banach vector bundle over $M$
admits a smooth linear connection.
\end{theorem}

\begin{proof}
Simultaneously refining an atlas of $M$ and a family of trivializations,
we obtain a cover $(U_\alpha)_{\alpha\in A}$ such that each $U_\alpha$ is
the domain of a chart $\varphi_\alpha$ and a trivialization
$\Phi_\alpha$. Choose a smooth partition of unity
$(\rho_\alpha)_{\alpha\in A}$ subordinate to it. On
$\boldsymbol{\mathcal{E}}\restriction_{U_\alpha}$ there is a local flat connection
$\nabla^\alpha$ whose coefficient in $(\varphi_\alpha,\Phi_\alpha)$ is zero.

Define the rule on the entire sheaf of sections. Let $W\subseteq M$
be open, let $\mathbf{X}\in\mathfrak X(W)$, and let
$\mathbf{s}\in\Gamma^\infty(\boldsymbol{\mathcal{E}}\restriction_W)$. On $W\cap U_\alpha$ set
\[
T_\alpha(\mathbf{X},\mathbf{s})
:=
\rho_\alpha\restriction_W\,
\nabla^\alpha_{\mathbf{X}\restriction_{W\cap U_\alpha}}
  (\mathbf{s}\restriction_{W\cap U_\alpha}).
\]
Since $\operatorname{supp}(\rho_\alpha)\subseteq U_\alpha$, this section
extends by zero to a smooth section of
$\boldsymbol{\mathcal{E}}\restriction_W$: every point of $W\setminus U_\alpha$ has a
neighborhood on which $\rho_\alpha$ vanishes identically. The family of these
extensions is locally finite. Therefore, the sum
\[
\nabla_{\mathbf{X}}\mathbf{s}:=\sum_{\alpha\in A}T_\alpha(\mathbf{X},\mathbf{s})
\]
is well defined and smooth on $W$. The same formula on an open set
$W_0\subseteq W$ is the restriction of the preceding one, so the rule is
compatible with restrictions.

Each $\nabla^\alpha$ is linear over functions in the direction and linear over
$\mathbb R$ in the section; multiplication by $\rho_\alpha$ and summation preserve
these two properties. For $f\in C^\infty(W)$, on each
$W\cap U_\alpha$ we have
\begin{align*}
\nabla_{\mathbf{X}}(f\mathbf{s})
&=
\sum_{\alpha\in A}\rho_\alpha
\bigl(\mathbf{X}(f)\mathbf{s}+f\nabla_{\mathbf{X}}^\alpha \mathbf{s}\bigr)\\
&=
\mathbf{X}(f)\mathbf{s}+f\nabla_{\mathbf{X}}\mathbf{s},
\end{align*}
since $\displaystyle\sum_{\alpha\in A}\rho_\alpha=1$. Finally, in any chart and
trivialization, each local flat connection transforms by
\eqref{eq:ley-transformacion-conexion-banach} into a smooth coefficient with
values in $\mathcal{L}(X,\mathcal{L}(F))$. The coefficient of $\nabla$ is the
locally finite sum of those coefficients multiplied by
$\rho_\alpha$; hence it is smooth in the operator norm. Thus
all the conditions of Definition~\ref{def:conexion-lineal-haz-banach} hold.
\end{proof}

The hypothesis of the preceding theorem is substantive in the Banach category.
Paracompactness alone does not produce smooth partitions of unity; the
hypotheses on the model space studied in
Section~\ref{sec:particiones-unidad-banach} are required.

\begin{proposition}[Induced connections]
\label{prop:conexiones-inducidas-haces-banach}
Let $M$ and $N$ be smooth Banach manifolds, and let $\boldsymbol{\mathcal E}\to M$ be a smooth Banach vector bundle.
\index{connection!dual on a Banach bundle}
\index{connection!pullback on a Banach bundle}
Let $\nabla$ be a connection on $\boldsymbol{\mathcal{E}}$. Then:
\begin{enumerate}[label=\textup{(\roman*)}]
\item the dual bundle $\boldsymbol{\mathcal{E}}'$ has the dual connection $\nabla'$ determined
by the identity
\begin{equation}
\label{eq:conexion-dual-banach}
(\nabla'_{\mathbf{X}}\boldsymbol{\lambda})(\mathbf{s})
=
\mathbf{X}\bigl(\boldsymbol{\lambda}(\mathbf{s})\bigr)-\boldsymbol{\lambda}(\nabla_{\mathbf{X}}\mathbf{s});
\end{equation}
here $U\subseteq M$ is open, $\mathbf{X}\in\mathfrak X(U)$,
$\boldsymbol{\lambda}\in\Gamma^\infty(\boldsymbol{\mathcal{E}}'\restriction_U)$, and
$\mathbf{s}\in\Gamma^\infty(\boldsymbol{\mathcal{E}}\restriction_U)$;
\item if $\boldsymbol{\mathcal{H}}$ has a connection $\nabla^{\boldsymbol{\mathcal{H}}}$, the bundle
$\boldsymbol{\mathcal{E}}\oplus\boldsymbol{\mathcal{H}}$ has the direct sum connection
\[
\nabla_{\mathbf{X}}^{\oplus}(\mathbf{s},\mathbf{t})
=
(\nabla_{\mathbf{X}}\mathbf{s},\nabla_{\mathbf{X}}^{\boldsymbol{\mathcal{H}}}\mathbf{t});
\]
\item if $f\colon N\longrightarrow M$ is smooth, the bundle $f^*\boldsymbol{\mathcal{E}}$ has a connection
$f^*\nabla$. Let $(V,\psi)$ be a chart of $N$ and $(U,\varphi)$ a chart of
$M$ such that $f(V)\subseteq U$, and let
$\Phi\colon \boldsymbol{\mathcal{E}}\restriction_U\longrightarrow U\times F$ be a
trivialization. Let $\widehat f:=\varphi\circ f\circ\psi^{-1}$, and denote
by $\Gamma$ the coefficient of $\nabla$ in $(\varphi,\Phi)$. Let $\widetilde\Phi$ denote the induced trivialization over $V$.
Let $\mathbf{Z}\in\mathfrak X(V)$ and
$\boldsymbol{\sigma}\in\Gamma^\infty((f^*\boldsymbol{\mathcal E})\restriction_V)$,
and let $Z$ and $\sigma$ be their local representations in the chart
$(V,\psi)$ and the trivialization $\widetilde\Phi$, respectively.
For $y\in\psi(V)$, the covariant derivative has the expression
\begin{equation}
\label{eq:conexion-pullback-banach}
\begin{aligned}
&\bigl((f^*\nabla)_{\mathbf{Z}}\boldsymbol{\sigma}\bigr)_{\psi,\widetilde\Phi}(y)\\
&\qquad=
D\sigma(y)[Z(y)]
+\Gamma(\widehat f(y))
\bigl(D\widehat f(y)[Z(y)]\bigr)\sigma(y).
\end{aligned}
\end{equation}
\end{enumerate}
\end{proposition}

\begin{proof}
For \eqref{eq:conexion-dual-banach} to define a covector, we must first
check its dependence on the test section. If
$h\in C^\infty(U)$, then
\begin{align*}
\mathbf{X}\bigl(\boldsymbol{\lambda}(h\mathbf{s})\bigr)-\boldsymbol{\lambda}(\nabla_{\mathbf{X}}(h\mathbf{s}))
&=\mathbf{X}\bigl(h\boldsymbol{\lambda}(\mathbf{s})\bigr)
  -\boldsymbol{\lambda}\bigl(\mathbf{X}(h)\mathbf{s}+h\nabla_{\mathbf{X}}\mathbf{s}\bigr)\\
&=h\bigl(\mathbf{X}(\boldsymbol{\lambda}(\mathbf{s}))-\boldsymbol{\lambda}(\nabla_{\mathbf{X}}\mathbf{s})\bigr).
\end{align*}
The right-hand side of \eqref{eq:conexion-dual-banach} is therefore
$C^\infty(U)$-linear in $\mathbf{s}$. To check also that its value
at a point depends only on the fiber, write $\ell$ and $s$ for the
representations of $\boldsymbol\lambda$ and $\mathbf s$, and $u=X_\varphi(x)$.
The derivative of evaluation, which is continuous and bilinear, gives
\[
 D\bigl(\ell(\cdot)(s(\cdot))\bigr)(x)[u]
 =D\ell(x)[u](s(x))+\ell(x)(Ds(x)[u]).
\]
Subtracting $\ell(x)(Ds(x)[u]+\Gamma(x)(u)s(x))$ cancels the terms involving $Ds$.
The result depends only on $s(x)$ and is given by the
continuous functional
\[
D\ell(x)[X_\varphi(x)]-
  \bigl(\Gamma(x)(X_\varphi(x))\bigr)'\ell(x).
\]
In particular, it is continuous in the fiber variable and depends smoothly on the
point. The coefficient of the dual connection is
\[
\Gamma'(x)(v)=-\bigl(\Gamma(x)(v)\bigr)',
\]
which is smooth in the operator norm because transposition
$\mathcal{L}(F)\longrightarrow\mathcal{L}(F')$ is continuous and linear.
Differentiating the identity
$\langle(g^{-1})'\ell,gw\rangle=\langle\ell,w\rangle$ and using
\eqref{eq:ley-transformacion-conexion-banach} gives exactly the transformation
law of $\Gamma'$ for the dual cocycle $(g^{-1})'$. Linearity
and the Leibniz rule in $\boldsymbol{\lambda}$ follow directly from
\eqref{eq:conexion-dual-banach}.

The direct sum connection satisfies the identities componentwise; in
simultaneous trivializations its coefficient is the diagonal operator
$\Gamma^{\boldsymbol{\mathcal{E}}}\oplus\Gamma^{\boldsymbol{\mathcal{H}}}$.

For the pullback, suppose that $N$ is modeled on a Banach space
$Y$. In the coordinates and induced trivialization of the statement, define
\[
\Gamma^f(y)(w)
:=
\Gamma(\widehat f(y))\bigl(D\widehat f(y)[w]\bigr),
\qquad w\in Y.
\]
Composition and evaluation of operators are smooth, so
$\Gamma^f\colon \psi(V)\longrightarrow\mathcal{L}(Y,\mathcal{L}(F))$ is
smooth in the operator norm. If the domain chart, the codomain
chart, and the trivialization are changed, substitute
\eqref{eq:ley-transformacion-conexion-banach} into this formula and use
\[
D(g\circ\widehat f)(y)[w]
=Dg(\widehat f(y))[D\widehat f(y)[w]].
\]
The result is \eqref{eq:ley-transformacion-conexion-banach} for the induced transition
map $g\circ\widehat f$; hence the local expressions
glue together. The resulting local formula is
\eqref{eq:conexion-pullback-banach}. Its linearity over functions in the
vector direction and
the Leibniz rule in $\sigma$ follow, respectively, from the linearity of
$D\widehat f(y)$ and the product rule for $D(h\sigma)$.
\end{proof}

The dual connection and the induced connection extend the constructions of
Proposition~\ref{conexion en tensores, duales y endomorfismos}. In finite dimensions the
coefficients are matrices; here the same proof additionally requires that
transposition, composition, and the coefficients depend continuously in the
operator norm.

\section{Differentiation along curves and parallel transport}
\label{sec:transporte-paralelo-banach}

Let $\gamma\colon I\longrightarrow M$ be a curve of class $C^1$. A field of class $C^1$ along
$\gamma$ is a section $\mathbf{V}\in\Gamma^1(\gamma^*\boldsymbol{\mathcal E})$. If
$\gamma(J)$ is contained in the domain of a chart and a trivialization,
write $x(t)$ for the representation of $\gamma$ and $v(t)$ for that of
$\mathbf{V}$. The local formula \eqref{eq:conexion-local-haz-banach} shows that
$(\nabla_{\mathbf{X}}\mathbf{s})(p)$ depends on $\mathbf{X}$ only through $\mathbf{X}(p)$. Hence, if
$u\in T_pM$, we write $\nabla_u\mathbf{s}$ for the common value of
$(\nabla_{\mathbf{X}}\mathbf{s})(p)$ among all local fields $\mathbf{X}$ such that $\mathbf{X}(p)=u$.

\begin{definition}[Covariant derivative along a curve]
\label{def:derivada-covariante-curva-banach}
Let $M$ be a smooth Banach manifold, let $\boldsymbol{\mathcal E}\to M$ be a smooth Banach vector bundle with connection, and let $\gamma\colon I\to M$ be a curve of class $C^1$.
\index{covariant derivative!along a curve on a Banach manifold}
The \textbf{covariant derivative of $\mathbf{V}$ along $\gamma$} is the field
$\frac{D^\nabla \mathbf{V}}{dt}$ whose representation in the trivialization $\Phi$ is defined
locally as follows. For $p$ in the domain of $\Phi$, denote by
\[
\Phi_p:=\operatorname{pr}_2\circ
\bigl(\Phi\restriction_{\boldsymbol{\mathcal E}_p}\bigr)\colon \boldsymbol{\mathcal E}_p\longrightarrow F
\]
the fiber coordinate. Then
\begin{equation}
\label{eq:derivada-covariante-curva-banach}
\Phi_{\gamma(t)}
\left(\left(\frac{D^\nabla \mathbf{V}}{dt}\right)(t)\right)
=
v'(t)+\Gamma(x(t))(x'(t))v(t).
\end{equation}
For $\gamma$ and $\mathbf V$ of class $C^1$, the covariant derivative is a
continuous field along $\gamma$. The field $\mathbf{V}$ is
\textbf{parallel} if $\frac{D^\nabla \mathbf{V}}{dt}=0$.
\end{definition}

\begin{proposition}[Invariance of differentiation along curves]
\label{prop:invariancia-derivada-curva-banach}
Let $M$ be a smooth Banach manifold, let $\boldsymbol{\mathcal E}\to M$ be a smooth Banach vector bundle with connection, and let $\gamma\colon I\to M$ be a curve of class $C^1$.
Expression \eqref{eq:derivada-covariante-curva-banach} does not depend on the
chart or the trivialization. It is linear in $\mathbf{V}$ and satisfies
\[
\frac{D^\nabla(a\mathbf{V})}{dt}
=
a'\mathbf{V}+a\frac{D^\nabla \mathbf{V}}{dt}
\]
for every function $a\in C^1(I)$. If $\mathbf{V}(t)=\mathbf{s}(\gamma(t))$ for a smooth local
section $\mathbf{s}$, then
\[
\frac{D^\nabla \mathbf{V}}{dt}(t)
=
\nabla_{\dot\gamma(t)}\mathbf{s}.
\]
\end{proposition}

\begin{proof}
Under a change of trivialization we have
$\widetilde v(t)=g(x(t))v(t)$. Therefore,
\[
\widetilde v'(t)
=
Dg(x(t))[x'(t)]v(t)+g(x(t))v'(t).
\]
Adding the connection term and using
\eqref{eq:ley-transformacion-conexion-banach} cancels the two terms containing
$Dg$, and the result is $g(x(t))$ times the right-hand side of
\eqref{eq:derivada-covariante-curva-banach}. Invariance under changes of
base coordinates is included in the same law. The remaining
assertions follow from the product rule and
\eqref{eq:conexion-local-haz-banach}.
\end{proof}

\begin{theorem}[Existence and uniqueness of parallel transport]
\label{teo:transporte-paralelo-haz-banach}
Let $M$ be a smooth Banach manifold, let $\boldsymbol{\mathcal E}\to M$ be a smooth Banach vector bundle with connection, and let $\gamma\colon I\to M$ be a curve of class $C^1$.
\index{parallel transport!in a Banach bundle}
Let $t_0\in I$ and $e_0\in\boldsymbol{\mathcal{E}}_{\gamma(t_0)}$. There exists a unique
parallel field $\mathbf{V}$ along $\gamma$ such that $\mathbf{V}(t_0)=e_0$. For
$t_0,t_1\in I$, the map
\[
P^\nabla_{\gamma,t_0\to t_1}\colon \boldsymbol{\mathcal{E}}_{\gamma(t_0)}
\longrightarrow\boldsymbol{\mathcal{E}}_{\gamma(t_1)},
\qquad
e_0\longmapsto \mathbf{V}(t_1),
\]
is a topological linear isomorphism. Moreover,
\begin{align*}
P^\nabla_{\gamma,t_1\to t_2}
\circ P^\nabla_{\gamma,t_0\to t_1}
&=P^\nabla_{\gamma,t_0\to t_2},\\
\bigl(P^\nabla_{\gamma,t_0\to t_1}\bigr)^{-1}
&=P^\nabla_{\gamma,t_1\to t_0}.
\end{align*}
for any $t_0,t_1,t_2\in I$.
\end{theorem}

\begin{proof}
In a trivialization over a subinterval $J\subseteq I$, the parallelism
condition is equivalent to the linear ODE in $F$
\begin{equation}
\label{eq:edo-transporte-paralelo-banach}
v'(t)=-A(t)v(t),
\qquad
A(t):=\Gamma(x(t))(x'(t))\in\mathcal{L}(F).
\end{equation}
The coefficient $A$ is continuous in the operator norm. The existence
and uniqueness theorem for ODEs in Banach spaces,
Theorem~\ref{teo:picard-lindelof-banach}, gives a unique local
solution. On every compact subinterval contained in $J$, the function
$\|A(t)\|$ is bounded, and the integral equation together with
Grönwall's lemma~\ref{lema: gronwall} gives
\[
\|v(t)\|_F
\leq
\|v(t_0)\|_F
\exp\left(\int_{t_0}^{t}\|A(\tau)\|\,d\tau\right)
\]
when $t\geq t_0$, with the analogous estimate in the reverse direction. Therefore,
the solution remains bounded on every compact subinterval. If a finite
endpoint $c$ of the maximal interval belonged to the interior of $J$, then,
on a compact set containing $t_0$ and $c$,
$A$ and $v$ would also be bounded. The equation would imply that $v'=-Av$ is bounded; hence $v$
would be uniformly Lipschitz and would have a limit in $F$ as $t$ tends to
$c$, because $F$ is complete. Applying the local theorem again with that limit
as initial value would extend the solution beyond $c$, a contradiction.
Thus, the solution exists throughout $J$. Covering the curve segment
between two times by finitely many trivializations and
using uniqueness on the overlaps gives a solution on all of
$I$.

The equation is linear in the initial value, and the preceding estimate, combined
over finitely many trivializations, proves continuity of the transport
operator. Uniqueness shows that transport from $t_0$ to $t_1$ followed by
transport from $t_1$ to $t_2$ produces the same parallel section as direct
transport. Taking $t_2=t_0$ gives the formula for the inverse.
\end{proof}

The notion of a section along a curve and the differentiation formula are
those of Definition
\ref{def:variedades-riemannianas-secciones-a-lo-largo-de-una-curva} and
Theorem~\ref{teo: derivada covariante a lo largo de una curva}. Likewise, the
composition and inversion laws correspond to
Theorem~\ref{teo: existencia y unicidad transporte paralelo} and
Proposition~\ref{prop: transporte paralelo isomorfismo}. The new feature here is
that the linear equation takes values in a Banach fiber; global existence
on the interval follows from the estimate in
Grönwall's lemma~\ref{lema: gronwall} and the
completeness of the fiber, rather than from finite dimensionality.

\begin{definition}[Curvature of a connection]
\label{def:curvatura-conexion-banach}
Let $M$ be a smooth Banach manifold, and let $\boldsymbol{\mathcal E}\to M$ be a smooth Banach vector bundle with connection $\nabla$.
\index{curvature!of a connection on a Banach bundle}
The \textbf{curvature} of $\nabla$ is the operator
\[
\mathbf{R}^\nabla(\mathbf{X},\mathbf{Y})\mathbf{s}
:=
\nabla_{\mathbf{X}}\nabla_{\mathbf{Y}}\mathbf{s}-
\nabla_{\mathbf{Y}}\nabla_{\mathbf{X}}\mathbf{s}-
\nabla_{[\mathbf{X},\mathbf{Y}]}\mathbf{s}.
\]
The connection is \textbf{flat} if $\mathbf{R}^\nabla=0$.
\end{definition}

\begin{proposition}[Tensoriality and local expression of curvature]
\label{prop:curvatura-local-conexion-banach}
Let $M$ be a smooth Banach manifold, and let $\boldsymbol{\mathcal E}\to M$ be a smooth Banach vector bundle with connection $\nabla$.
Curvature is $C^\infty(M)$-linear in $\mathbf{X}$, $\mathbf{Y}$, and $\mathbf{s}$, and alternating in
$\mathbf{X},\mathbf{Y}$. The local expression obtained below
also shows that its value at $p$ depends only on
$\mathbf X(p)$, $\mathbf Y(p)$, and $\mathbf s(p)$. Thus, for each $p\in M$
it determines a continuous trilinear map
\[
\mathbf{R}_p^\nabla\colon T_pM\times T_pM\times\boldsymbol{\mathcal E}_p\longrightarrow\boldsymbol{\mathcal E}_p.
\]
Equivalently,
\[
\mathbf{R}_p^\nabla
\in
\mathcal L^2_{\mathrm{alt}}
\bigl(T_pM;\mathcal L(\boldsymbol{\mathcal E}_p)\bigr).
\]
In a chart and a trivialization, if $x$ represents $p$, then
\begin{equation}
\label{eq:curvatura-local-conexion-banach}
\mathbf{R}_x^\nabla(u,v)w
=
\Omega(x)(u,v)w,
\end{equation}
where
\begin{align*}
\Omega(x)(u,v)
={}&D\Gamma(x)[u](v)-D\Gamma(x)[v](u)\\
&+\Gamma(x)(u)\Gamma(x)(v)
-\Gamma(x)(v)\Gamma(x)(u).
\end{align*}
Under a change of chart and trivialization as in
Proposition~\ref{prop:transformacion-conexion-haz-banach}, we have
\begin{equation}
\label{eq:transformacion-curvatura-banach}
\widetilde\Omega(\chi(x))(D\chi(x)u,D\chi(x)v)
=
g(x)\Omega(x)(u,v)g(x)^{-1}.
\end{equation}
\end{proposition}

\begin{proof}
Antisymmetry is immediate. For $f\in C^\infty(M)$, we use
\[
[f\mathbf{X},\mathbf{Y}]=f[\mathbf{X},\mathbf{Y}]-\mathbf{Y}(f)\mathbf{X}
\]
and the Leibniz rules for the connection to obtain
$\mathbf{R}^\nabla(f\mathbf{X},\mathbf{Y})\mathbf{s}=f\mathbf{R}^\nabla(\mathbf{X},\mathbf{Y})\mathbf{s}$. Linearity in the second direction
follows from antisymmetry. Likewise, on expanding
$\mathbf{R}^\nabla(\mathbf{X},\mathbf{Y})(f\mathbf{s})$, the terms involving first derivatives of $f$ cancel, and
the remaining second-order term is
\[
\bigl(\mathbf{X}(\mathbf{Y}(f))-\mathbf{Y}(\mathbf{X}(f))-[\mathbf{X},\mathbf{Y}](f)\bigr)\mathbf{s}=0.
\]
These identities prove linearity over functions. Pointwise dependence
and continuity are checked directly using the coefficients of
the connection. Denote by $U,V$ and $s$ the local representations of the
two fields and the section; in the following equalities all maps and
their derivatives are evaluated at $x$. Applying
\eqref{eq:conexion-local-haz-banach} twice gives
\begin{align*}
 (\nabla_{\mathbf X}\nabla_{\mathbf Y}\mathbf s)_{\varphi,\Phi}
 ={}&D^2s[U,V]+Ds[DV[U]]+D\Gamma[U](V)s+\Gamma(DV[U])s\\
 &+\Gamma(V)Ds[U]+\Gamma(U)Ds[V]+\Gamma(U)\Gamma(V)s.
\end{align*}
The expression for $\nabla_{\mathbf Y}\nabla_{\mathbf X}\mathbf s$ is obtained by
interchanging $U$ and $V$. Subtracting them cancels the terms involving $D^2s$
by symmetry, and the two terms in which a derivative of $s$ is
multiplied by $\Gamma$ appear on both sides with the same sign.
The terms containing $DU$ or $DV$ sum to
\[
 Ds[DV[U]-DU[V]]+\Gamma(DV[U]-DU[V])s
 =(\nabla_{[\mathbf X,\mathbf Y]}\mathbf s)_{\varphi,\Phi}.
\]
Subtracting this last expression leaves
\eqref{eq:curvatura-local-conexion-banach}, now obtained for arbitrary fields and
sections. They may therefore be replaced by constant fields
with the same values without changing the result. For $u,v\in X$ and $w\in F$,
the formula gives the bound
\[
 \|\Omega(x)(u,v)w\|_F
 \leq\bigl(2\|D\Gamma(x)\|+2\|\Gamma(x)\|^2\bigr)
       \|u\|_X\|v\|_X\|w\|_F.
\]
The norms of $\Gamma(x)$ and $D\Gamma(x)$ are those of
$\mathcal L(X,\mathcal L(F))$ and
$\mathcal L(X,\mathcal L(X,\mathcal L(F)))$, respectively.
This proves trilinear continuity; smoothness in $x$ follows from the
smoothness of $\Gamma$ and continuous bilinear composition of operators. Finally, since
$\mathbf{R}^\nabla(\mathbf{X},\mathbf{Y})\mathbf{s}$ is a section and is tensorial in its three arguments, its
local representations satisfy
\eqref{eq:transformacion-curvatura-banach}; equivalently, the identity
follows by substituting \eqref{eq:ley-transformacion-conexion-banach} and
canceling the terms containing derivatives of $g$ and $\chi$.
\end{proof}

The definition and proof of tensoriality correspond to
Definition~\ref{def:variedades-riemannianas-variedad-5} and
Proposition~\ref{prop:variedades-riemannianas-tensor-curvatura-riemann-define-campo-tensorial}.
The algebraic cancellations are exactly the same. What must
additionally be checked in the Banach setting is that
$D\Gamma(x)[u](v)$ and the operator products in
\eqref{eq:curvatura-local-conexion-banach} depend smoothly on $x$ in the
operator norm; this follows from the smoothness of $\Gamma$ and the
continuous bilinearity of evaluation and composition in
$\mathcal L(F)$.

\chapter{Hilbert geometry and Finsler structures}
\label{cap:geometria-hilbert-finsler}

The existence of a Banach tangent space does not by itself provide a
notion of length. It is useful to distinguish three structures that are often
confused: a strong Riemannian metric on a Hilbert manifold, a
weak Riemannian metric, and a Banach--Finsler structure. The first
allows classical Riemannian geometry to be transported locally through the
Riesz isomorphism; the second may lack even gradients, a
Levi--Civita connection, or a distance that separates points; the third retains
the local metric control needed for completeness and deformation
flows. At the end we distinguish this theory from classical Finsler geometry,
whose vertical convexity leads to the geodesic spray and the Chern connection.
Manifolds are understood according to the convention in
Definition~\ref{def:variedad-banach}. In Riemannian geometry and classical
Finsler geometry we work with smooth differentiable
structures; the lower regularity needed for variational applications
will be specified in the Banach--Finsler section.

\begin{semblanzaHistorica}{What remains of Riemannian geometry in infinite dimensions}
In finite dimensions, every Riemannian metric identifies vectors and covectors and produces a distance locally equivalent to the coordinate distance. On a Hilbert manifold, that identification may fail if the metric is weak. This difference, studied in modern infinite-dimensional geometry, explains why some gradients do not exist, why the distance can degenerate, and why the Hopf--Rinow theorem loses its usual equivalences. Strong and Banach--Finsler structures recover part of the necessary control, but each conclusion must be proved under its own hypotheses.
\end{semblanzaHistorica}

\section{Hilbert manifolds and strong or weak metrics}
\label{sec:metricas-hilbert-fuertes-debiles}

\begin{definition}\label{def:variedad-hilbert}\glsadd{variedad-hilbert}\glsadd{metrica-riemanniana-debil}\glsadd{metrica-riemanniana-fuerte}
\index{Hilbert manifold}\index{Riemannian metric!strong}\index{Riemannian metric!weak}
A \textbf{Hilbert manifold} is a Banach manifold whose model space
is a real Hilbert space $H$. Let $M$ be a Hilbert manifold modeled
on $H$. A \textbf{smooth weak Riemannian metric} on $M$ is a family
$\mathbf{g}=(\mathbf{g}_x)_{x\in M}$ of continuous positive definite inner products on
$T_xM$ such that, in every trivialization of $TM$, the map
\[
 x\longmapsto \mathbf{g}_x\in\mathcal L_s^2(H;\mathbb R)
\]
is smooth. Its musical homomorphism is
\[
 \flat_{\mathbf{g}}\colon TM\longrightarrow T'M,
 \qquad \flat_{\mathbf{g}}(v)=\mathbf{g}(v,\,\cdot\,).
\]
Here $T'M$ is the bundle of continuous duals from
Definition~\ref{def:haz-cotangente-banach}; in coordinates its fibers are
modeled on $H'$.
The metric is \textbf{strong} if each norm $|v|_{\mathbf{g},x}=\mathbf{g}_x(v,v)^{\frac{1}{2}}$
induces the given topology of $T_xM$; it is \textbf{strictly weak} if it is not
strong. Thus, ``weak'' here means ``not necessarily strong.''
Positivity implies that $\flat_{\mathbf{g},x}$ is injective, but not that it is
surjective.
\end{definition}

Fix an inner product $\langle\,\cdot\,,\cdot\,\rangle_H$ on the model
space of a Hilbert manifold. In a chart there is a unique
positive, injective, self-adjoint operator $G(x)\in\mathcal L(H)$ such that
\begin{equation}\label{eq:operador-metrico-hilbert}
 \mathbf{g}_x(u,v)=\langle G(x)u,v\rangle_H.
\end{equation}
This is the representation in
Proposition~\ref{prop:formas-bilineales-simetricas-hilbert}, applied at each
point.
The map $x\mapsto G(x)$ is smooth. After identifying $H'$ with $H$ by
Riesz, $\flat_{\mathbf{g},x}$ is represented by $G(x)$. In particular, its image
is dense. Indeed, since $G(x)$ is self-adjoint,
$(\operatorname{Ran}G(x))^\perp=\ker(G(x)^*)=\ker G(x)=\{0\}$ and hence
\[
 \overline{\operatorname{Ran}G(x)}
 =\bigl(\ker G(x)\bigr)^\perp=H.
\]
The range may nevertheless be proper and nonclosed.

\begin{proposition}\label{prop:caracterizaciones-metrica-fuerte}
For a weak Riemannian metric $\mathbf{g}$ on a Hilbert manifold, the following are
equivalent:
\begin{enumerate}[label=\textup{(\roman*)}]
 \item $\mathbf{g}$ is strong;
 \item $\flat_{\mathbf{g},x}\colon T_xM\longrightarrow T'_xM$ is a topological
       isomorphism for every $x\in M$;
 \item $\flat_{\mathbf{g}}\colon TM\longrightarrow T'M$ is a smooth vector bundle
       isomorphism;
 \item about each $x_0$ there exist a chart and constants $0<c\leq C$ such
       that
       \[
       c\|v\|_H\leq |v|_{\mathbf{g},x}\leq C\|v\|_H
       \quad\text{for every $x$ in a neighborhood of $x_0$ and $v\in H$.}
       \]
 \item about every point there is a Hilbert chart in which the
       operator $G(x)$ of \eqref{eq:operador-metrico-hilbert} belongs to
       $\operatorname{GL}(H)$ for every $x$ in its domain; this property
       then holds in every chart.
\end{enumerate}
Moreover, if a nonempty Banach manifold admits a strong metric, its model space
is isomorphic to a Hilbert space.
\end{proposition}

\begin{proof}
The norm of \eqref{eq:operador-metrico-hilbert} induces the topology of $H$ if
and only if there exist $c,C>0$ such that
\[
 c^2\|v\|_H^2\leq\langle G(x)v,v\rangle_H\leq C^2\|v\|_H^2.
\]
The lower inequality is equivalent to invertibility of $G(x)$. Indeed, if it
holds, $G(x)^{\frac{1}{2}}$ is bounded below and has closed dense image,
so it is surjective. The converse follows from
continuity of $G(x)^{-\frac{1}{2}}$. This proves the pointwise equivalences between
(i), (ii), and (v).

Since $\operatorname{GL}(H)$ is open and inversion is smooth, continuity
of $G$ and $G^{-1}$ gives, after shrinking the chart, uniform bounds
for $\|G(x)\|$ and $\|G(x)^{-1}\|$. Specifically, if these bounds are $A$ and $B$,
respectively, we may take $C=A^{\frac{1}{2}}$ and $c=B^{-\frac{1}{2}}$ in (iv). The same
observation shows that the fiberwise inverse
$\sharp_{\mathbf{g}}:=\flat_{\mathbf{g}}^{-1}$ is smooth. Thus, (iii) is equivalent to the preceding
conditions.

For the last assertion, denote the model space by $X$, fix
$x_0\in M$, and transport $\mathbf{g}_{x_0}$ to $X$. Since the given norm and the inner product norm
$|\,\cdot\,|_{\mathbf{g},x_0}$ induce the same topology, there exist $c_0,C_0>0$ such
that
\[
 c_0\|v\|_X\leq |v|_{\mathbf{g},x_0}\leq C_0\|v\|_X,
 \qquad v\in X.
\]
If $(v_n)$ is Cauchy in $|\,\cdot\,|_{\mathbf{g},x_0}$, the left inequality
makes it Cauchy in $\|\,\cdot\,\|_X$; it then converges to some $v\in X$, and
the right inequality turns this convergence into convergence in
$|\,\cdot\,|_{\mathbf{g},x_0}$. Consequently, $X$, equipped with this inner
product, is a Hilbert space, and the identity between the two norms is a
topological linear isomorphism.
\end{proof}

\begin{example}[Hierarchy of Sobolev metrics]
\label{ej:metricas-sobolev-fuerte-debil}
Let $\mathbf{E}\longrightarrow K$ be a real vector bundle of positive finite rank, equipped with a bundle metric and a compatible
connection, over a compact
Riemannian manifold without boundary $K$, and let
$A=I+(\nabla^{\mathbf{E}})^*\nabla^{\mathbf{E}}$. Suppose that $K$ is nonempty. For $s>0$ and
$0\leq r\leq s$, on the Hilbert space $H^s(K,\mathbf{E})$ consider the
constant metric
\begin{equation}\label{eq:metrica-sobolev-orden-r}
 \mathbf{g}^{(r)}_{\mathbf{u}}(\mathbf{h},\mathbf{k})
 :=\bigl\langle A^{\frac{r}{2}}\mathbf{h},A^{\frac{r}{2}}\mathbf{k}
 \bigr\rangle_{L^2(K,\mathbf{E})}.
\end{equation}
Here we use the spectral norm
$\|\mathbf h\|_{H_A^t(K,\mathbf{E})}
:=\|A^{\frac{t}{2}}\mathbf h\|_{L^2(K,\mathbf{E})}$. For each $t\geq0$ there exist constants
$c_t,C_t>0$, depending only on $K$, the metrics, the connection, and
$t$, such that
\[
 c_t\|\mathbf{h}\|_{H^t(K,\mathbf{E})}
 \leq\|\mathbf{h}\|_{H_A^t(K,\mathbf{E})}
 \leq C_t\|\mathbf{h}\|_{H^t(K,\mathbf{E})},
 \qquad \mathbf{h}\in H^t(K,\mathbf{E}).
\]
For even integer orders, these are the estimates of
Corollary~\ref{cor:potencias-laplaciano-bochner-dominios-sobolev}; the remaining
orders follow through the spectral scale and interpolation.
Specifically, for $0<t<2m$, with $m\in\mathbb N$ and parameter $\theta=t/(2m)$, interpolate between
$L^2(K,\mathbf E)$ and the domain of $A^m$.
Proposition~\ref{prop:compatibilidad-espectral-dual-interpolacion-bessel-haces}
identifies the interpolation space with the domain of $A^{t/2}$, whereas
Proposition~\ref{prop:dualidad-sobolev-haces-geometria-acotada}
identifies it with $H^t(K,\mathbf E)$. Compactness of $K$ allows the
bounded geometry hypotheses of that result to be applied. This is also the case
$p=2$ of Theorem~\ref{teo:bessel-calor-triebel-local-haces-geometria-acotada}.
Continuity of the embedding
$H^s(K,\mathbf{E})\hookrightarrow H^r(K,\mathbf{E})$
shows that $\mathbf{g}^{(r)}$ is a smooth weak metric. For $r=s$ it is strong, since
its norm is one of the norms defining $H^s(K,\mathbf{E})$.

Suppose that $K$ has positive dimension. If
$A\mathbf{e}_j=\lambda_j\mathbf{e}_j$, with $(\mathbf{e}_j)_{j\in\mathbb N}$ orthonormal in $L^2(K,\mathbf{E})$ and
$\lambda_j\to\infty$ as $j\to\infty$, then
\[
 \mathbf{h}_j:=\lambda_j^{-\frac{s}{2}}\mathbf{e}_j,
 \qquad \|\mathbf{h}_j\|_{H_A^s(K,\mathbf{E})}=1,
 \qquad |\mathbf{h}_j|_{\mathbf{g}^{(r)}}=\lambda_j^{\frac{r-s}{2}}.
\]
Consequently, $\mathbf{g}^{(r)}$ is strictly weak for $r<s$. With respect to the
inner product of $H^s(K,\mathbf{E})$, its metric operator is $A^{r-s}$, which is
injective and has dense image, but is not surjective. In particular, there are
continuous covectors on $H^s(K,\mathbf{E})$ that have no gradient with respect
to $\mathbf{g}^{(r)}$.
\end{example}

\begin{example}[The $L^2$ metric]
\label{ej:metrica-L2-debil}
The case $r=0$ of \eqref{eq:metrica-sobolev-orden-r} is
\[
 \mathbf{g}^{(0)}_{\mathbf{u}}(\mathbf{h},\mathbf{k})
 =\int_K\langle \mathbf{h},\mathbf{k}\rangle_{\mathbf{E}}\,d\lambda_{\mathbf{g}}.
\]
Its musical map is the inclusion
\[
 H^s(K,\mathbf{E})\longrightarrow H^{-s}(K,\mathbf{E}),
 \qquad \mathbf{h}\longmapsto
 \bigl[\mathbf{k}\longmapsto(\mathbf{h},\mathbf{k})_{L^2(K,\mathbf{E})}\bigr].
\]
If $\dim K>0$, the image is proper, as checked in the preceding
example for $r=0<s$. Thus, using an $L^2$ integral on a configuration
space of regularity $H^s(K,\mathbf{E})$ does not automatically make the metric
strong. In contrast, the $\mathbf{g}^{(s)}$ metric does identify
$H^s(K,\mathbf{E})$ with its dual through Riesz.
\end{example}

The obstruction from the musical map is independent of the notion of smoothness chosen for
function spaces. In particular, passing to convenient calculus on
Fréchet spaces, developed in \cite{KrieglMichor1997}, does not make
a musical map surjective if it was not already. Convenient theory and strong
Hilbert theory solve different problems. A modern comparison of
categories of infinite-dimensional calculus and weak Riemannian geometry
can be found in \cite[Chap.~4 and Appendices~C and E]{Schmeding2022}. In this
book, however, all maps between open subsets of Banach spaces are interpreted
within the Fréchet calculus established in
Chapter~\ref{cap:calculo-diferencial-banach}; convenient calculus
is not adopted as a premise.

\section{Connections, geodesics, and variation in the strong case}
\label{sec:conexion-geodesicas-hilbert}

Let $\mathbf{g}$ be strong, and write $Dg(x)[u](v,w)$ for the derivative of its
local representative $g$. For $x\in U$ and $u,v\in H$, define the covector
\begin{equation}\label{eq:covector-koszul-local-hilbert}
 \mathscr K_x(u,v)(w)
 :=\frac12\bigl(Dg(x)[u](v,w)+Dg(x)[v](u,w)-Dg(x)[w](u,v)\bigr).
\end{equation}
With the multilinear operator norm we have
\begin{equation}\label{eq:cota-koszul-local-hilbert}
 \|\mathscr K_x(u,v)\|_{H'}
 \leq\frac32
 \|Dg(x)\|_{\mathcal L(H;\mathcal L_s^2(H;\mathbb R))}
 \|u\|_H\|v\|_H.
\end{equation}
Indeed, each of the three summands in
\eqref{eq:covector-koszul-local-hilbert} is bounded by the product
appearing on the right before the factor $\frac{1}{2}$. Continuous regrouping of
multilinear operators then shows that
$x\mapsto\mathscr K_x\in\mathcal L_s^2(H;H')$ is smooth. In particular, the
map $(x,u,v)\mapsto\mathscr K_x(u,v)$ is smooth, bilinear, and symmetric in
$(u,v)$. Since $\flat_{\mathbf{g},x}$ is
surjective, there exists a unique vector
$\Gamma_x(u,v)=\sharp_{\mathbf{g},x}\mathscr K_x(u,v)$ such that
\begin{equation}\label{eq:christoffel-hilbert}
2g_x\bigl(\Gamma_x(u,v),w\bigr)
=Dg(x)[u](v,w)+Dg(x)[v](u,w)-Dg(x)[w](u,v).
\end{equation}

Definition~\ref{def conexion compatible y libre de torsion} and
Theorem~\ref{teo:variedades-riemannianas-teorema-fundamental-de-la-geometria-riemanniana}
are the finite-dimensional model for the following construction. Its local
expression and coordinate change law appear in
Corollary~\ref{expresion local christoffel levi civita} and
Lemma~\ref{lema: transformacion metrica Christoffel}. The Koszul identity and
the calculation of chart changes are the same; the new step is that
its right-hand side is initially an element of $H'$, and only the strength of the
metric allows $\sharp_{\mathbf{g}}$ to be applied to it in the operator norm.

\begin{theorem}[Strong Levi--Civita theorem]\label{teo:levi-civita-hilbert-fuerte}
Every smooth strong Riemannian metric on a Hilbert manifold has a
unique torsion-free connection $\nabla$ compatible with $\mathbf{g}$. If $X,Y$
are the local representatives of the fields $\mathbf{X},\mathbf{Y}$, then
\[
 (\nabla_{\mathbf{X}}\mathbf{Y})_\varphi(x)
 =DY(x)X(x)+\Gamma_x(X(x),Y(x)),
\]
where $\Gamma$ is determined by \eqref{eq:christoffel-hilbert}.
Equivalently, for smooth fields $\mathbf{X},\mathbf{Y},\mathbf{Z}$ it satisfies the Koszul formula
\begin{align*}
2\mathbf{g}(\nabla_{\mathbf{X}}\mathbf{Y},\mathbf{Z})
=\mathbf{X}\mathbf{g}(\mathbf{Y},\mathbf{Z})
 +\mathbf{Y}\mathbf{g}(\mathbf{X},\mathbf{Z})
 -\mathbf{Z}\mathbf{g}(\mathbf{X},\mathbf{Y})+\mathbf{g}([\mathbf{X},\mathbf{Y}],\mathbf{Z})
 -\mathbf{g}([\mathbf{Y},\mathbf{Z}],\mathbf{X})
 -\mathbf{g}([\mathbf{X},\mathbf{Z}],\mathbf{Y}).
\end{align*}
If $\tau$ is a coordinate change between open subsets of Hilbert spaces,
the corresponding symbols satisfy the full law
\begin{equation}\label{eq:ley-transformacion-christoffel-hilbert}
 \widetilde\Gamma_{\tau(x)}\bigl(D\tau(x)u,D\tau(x)v\bigr)
 =D\tau(x)\Gamma_x(u,v)-D^2\tau(x)(u,v).
\end{equation}
Equivalently, if $\sigma=\tau^{-1}$ and $y=\tau(x)$,
\begin{align*}
 \widetilde\Gamma_y(a,b)
 =D\tau(\sigma(y))\Gamma_{\sigma(y)}
       \bigl(D\sigma(y)a,D\sigma(y)b\bigr)+D\tau(\sigma(y))D^2\sigma(y)(a,b).
\end{align*}
\end{theorem}

\begin{proof}
In a chart, the musical map depends smoothly on the point and takes values
in the open set of
linear isomorphisms $H\longrightarrow H'$. Since inversion is smooth on that open set,
$x\mapsto\sharp_{\mathbf{g},x}\in\mathcal L(H';H)$ is smooth. The preceding construction
can be written as
\[
 \Gamma_x=\sharp_{\mathbf{g},x}\circ\mathscr K_x.
\]
Continuous composition of operators, together with
\eqref{eq:cota-koszul-local-hilbert}, shows that
$x\mapsto\Gamma_x\in\mathcal L_s^2(H;H)$ is smooth. The local formula explicitly satisfies,
for $f,h\in C^\infty(U)$ and $a,b\in\mathbb R$,
\begin{align*}
 \nabla_{f\mathbf{X}+h\mathbf{Z}}\mathbf{Y}&=f\nabla_{\mathbf{X}}\mathbf{Y}+h\nabla_{\mathbf{Z}}\mathbf{Y},\\
 \nabla_{\mathbf{X}}(a\mathbf{Y}+b\mathbf{Z})&=a\nabla_{\mathbf{X}}\mathbf{Y}+b\nabla_{\mathbf{X}}\mathbf{Z},\\
 \nabla_{\mathbf{X}}(f\mathbf{Y})&=\mathbf{X}(f)\mathbf{Y}+f\nabla_{\mathbf{X}}\mathbf{Y},
\end{align*}
since $DY$ is linear in the direction, $\Gamma_x$ is bilinear, and, if
$\widehat f=f\circ\varphi^{-1}$, the product rule gives
$D(\widehat fY)[X]=D\widehat f[X]Y+\widehat fDY[X]$. These are therefore exactly the axioms of
Definition~\ref{def:conexion-lineal-haz-banach}, once invariance of the
local expressions has been checked.

Now let the representatives of the fields in the second coordinates be
$\widetilde X(y)=D\tau(\sigma(y))X(\sigma(y))$ and
$\widetilde Y(y)=D\tau(\sigma(y))Y(\sigma(y))$. The chain rule gives
\begin{align*}
 D\widetilde Y(\tau(x))\widetilde X(\tau(x))
 ={}&D^2\tau(x)(X(x),Y(x))\\
 &+D\tau(x)DY(x)X(x).
\end{align*}
Therefore, the equality
$\widetilde\nabla_{\widetilde X}\widetilde Y
=D\tau(\nabla_XY)$ between local expressions is equivalent precisely to
\eqref{eq:ley-transformacion-christoffel-hilbert}. To verify it from the metric, fix $x$ and set
$y=\tau(x)$, $A=D\tau(x)$, and $B=D^2\tau(x)$. The metric transformation
is $g_x(v,w)=\widetilde g_y(Av,Aw)$. Differentiating in the direction $u$,
\begin{align*}
 Dg(x)[u](v,w)
 ={}&D\widetilde g(y)[Au](Av,Aw)\\
 &+\widetilde g_y(B(u,v),Aw)+\widetilde g_y(Av,B(u,w)).
\end{align*}
In the combination defining $2\mathscr K_x(u,v)(w)$, the terms involving
$B(u,w)$ cancel those involving $B(w,u)$, and those involving $B(v,w)$ cancel those involving
$B(w,v)$. The two remaining terms are equal by symmetry of $B$.
Hence,
\[
 \mathscr K_x(u,v)(w)
 =\widetilde{\mathscr K}_y(Au,Av)(Aw)
  +\widetilde g_y(B(u,v),Aw).
\]
Substituting $\mathscr K_x(u,v)(w)=g_x(\Gamma_x(u,v),w)$ and the analogous
equality for $\widetilde{\mathscr K}$ gives
\[
 \widetilde g_y\bigl(A\Gamma_x(u,v)-\widetilde\Gamma_y(Au,Av)-B(u,v),Aw\bigr)=0.
\]
Since $A$ is surjective and $\widetilde g_y$ is nondegenerate, the first
argument must be zero. This is
\eqref{eq:ley-transformacion-christoffel-hilbert}; thus the local
expressions glue together and define a connection.

Symmetry of $\Gamma_x$ shows, in coordinates, that
\[
 (\nabla_{\mathbf{X}}\mathbf{Y}-\nabla_{\mathbf{Y}}\mathbf{X})_\varphi
 =DY[X]-DX[Y]=[\mathbf{X},\mathbf{Y}]_\varphi,
\]
so the torsion vanishes. On the other hand, the product rule gives
\[
 (\mathbf{X}\mathbf{g}(\mathbf{Y},\mathbf{Z}))\circ\varphi^{-1}=Dg[X](Y,Z)
 +g(DY[X],Z)
 +g(Y,DZ[X]).
\]
Adding the identities \eqref{eq:christoffel-hilbert} for
$(X,Y,Z)$ and $(X,Z,Y)$, the terms $Dg[Y](X,Z)$ and $Dg[Z](X,Y)$
appear once with each sign; by symmetry of $g$ they cancel, leaving
\[
 Dg[X](Y,Z)
 =g(\Gamma(X,Y),Z)
 +g(Y,\Gamma(X,Z)).
\]
The last two formulas prove metric compatibility. Conversely,
compatibility and the absence of torsion imply the Koszul formula through
the same algebraic calculation as in
Theorem~\ref{teo:variedades-riemannianas-teorema-fundamental-de-la-geometria-riemanniana}.
If $\nabla$ and $\widetilde\nabla$ had the same properties, the Koszul
formula would give
$\mathbf{g}((\nabla_{\mathbf{X}}-\widetilde\nabla_{\mathbf{X}})\mathbf{Y},\mathbf{Z})=0$
for every field $\mathbf{Z}$; nondegeneracy of $\mathbf{g}$ establishes uniqueness.
\end{proof}

\begin{proposition}[Musical criterion in the weak case]
\label{prop:criterio-levi-civita-debil}
Let $\mathbf{g}$ be a weak Riemannian metric on a Hilbert manifold. A
Levi--Civita connection for $\mathbf{g}$ exists if and only if, in an atlas, the covectors in
\eqref{eq:covector-koszul-local-hilbert} belong to
$\operatorname{Im}\flat_{\mathbf{g},x}$ and their unique preimages
\[
 \Gamma_x(u,v):=\flat_{\mathbf{g},x}^{-1}\mathscr K_x(u,v)
\]
depend smoothly on $(x,u,v)$. Equivalently, smoothness of
$x\mapsto\Gamma_x$ with values in $\mathcal L_s^2(H;H)$ equipped with its operator
norm is required. When the connection exists, it is unique.
\end{proposition}

\begin{proof}
Necessity follows by applying the Koszul formula to constant fields in
a chart. We check the equivalence of the two formulations of
smoothness. A family smooth in norm produces smooth joint
evaluation through continuous multilinear evaluation. Conversely, if
$H_0(x,u,v):=\Gamma_x(u,v)$ is smooth and bilinear in $(u,v)$, then
\[
 \Gamma_x(u,v)=D_vD_uH_0(x,0,0)[u,v].
\]
The Fréchet derivative on the right depends smoothly on $x$ in the
bilinear operator norm. Hence $\Gamma$ has exactly the
regularity required in Definition~\ref{def:conexion-lineal-haz-banach}.
The hypothesis therefore allows coefficients of a connection to be defined. Under a change of
coordinates, both sides of
\eqref{eq:ley-transformacion-christoffel-hilbert} agree, after application of
$\flat_{\mathbf{g}}$, by the chain rule. Since $\flat_{\mathbf{g}}$ is injective,
this law holds and the local expressions define a connection. The proof
of Theorem~\ref{teo:levi-civita-hilbert-fuerte} establishes compatibility,
absence of torsion, and uniqueness from the Koszul formula and
injectivity of the musical map.
\end{proof}

The criterion explains exactly the difficulty in the weak case. The constant metric
$\mathbf{g}^{(r)}$ of Example~\ref{ej:metricas-sobolev-fuerte-debil} has
$D\mathbf{g}^{(r)}=0$ and therefore admits the flat connection $\Gamma=0$, although it is
strictly weak when $r<s$. For a varying weak metric, however,
$\mathscr K_x(u,v)$ may lie outside the image of the musical map. There is then
no vector that can play the role of $\Gamma_x(u,v)$.

Let $\gamma\colon [a,b]\longrightarrow M$ be a $C^1$ curve, and let $\mathbf{V}$ be a
$C^1$ field along $\gamma$. Apply
Definition~\ref{def:derivada-covariante-curva-banach} to the bundle $TM$ and write
$\nabla_t\mathbf{V}:=\frac{D^\nabla \mathbf{V}}{dt}$. In a chart,
\[
 \nabla_t\mathbf{V}:=\dot V(t)+\Gamma_{\gamma(t)}(\dot\gamma(t),V(t)),
\]
Independence of the chart and the Leibniz rule are given by
Proposition~\ref{prop:invariancia-derivada-curva-banach}. Metric compatibility
adds, for every field $\mathbf{W}$ along $\gamma$, the second
of the rules
\[
 \nabla_t(f\mathbf{V})=f'\mathbf{V}+f\nabla_t\mathbf{V},
 \qquad
 \frac d{dt}\mathbf{g}(\mathbf{V},\mathbf{W})=\mathbf{g}(\nabla_t\mathbf{V},\mathbf{W})+\mathbf{g}(\mathbf{V},\nabla_t\mathbf{W}).
\]
Theorem~\ref{teo:transporte-paralelo-haz-banach}, applied to $TM$, gives
existence and uniqueness of fields satisfying $\nabla_t\mathbf{V}=0$, as well as the
topological linear isomorphisms
\[
 P_{s,t}^{\gamma}:=P^\nabla_{\gamma,s\to t}\colon
 T_{\gamma(s)}M\longrightarrow T_{\gamma(t)}M,
 \qquad
 P_{r,t}^{\gamma}P_{s,r}^{\gamma}=P_{s,t}^{\gamma},
 \qquad
 (P_{s,t}^{\gamma})^{-1}=P_{t,s}^{\gamma}.
\]
To identify the additional property supplied by Levi--Civita, let $\mathbf{U}$ and $\mathbf{V}$ be the parallel fields with initial values
$u$ and $v$ at time $s$. Metric compatibility gives
\[
 \frac d{dt}\mathbf{g}_{\gamma(t)}(\mathbf{U}(t),\mathbf{V}(t))
 =\mathbf{g}(\nabla_t\mathbf{U},\mathbf{V})+\mathbf{g}(\mathbf{U},\nabla_t\mathbf{V})=0.
\]
Consequently, for $s,t\in[a,b]$ and $u,v\in T_{\gamma(s)}M$,
$\mathbf{g}_{\gamma(t)}(P_{s,t}^{\gamma}u,P_{s,t}^{\gamma}v)
=\mathbf{g}_{\gamma(s)}(u,v)$ and parallel transport is an isometry. This is the
same deduction that, in finite dimensions, completes
Theorems~\ref{teo: derivada covariante a lo largo de una curva}
and~\ref{teo: existencia y unicidad transporte paralelo} and
Proposition~\ref{prop: transporte paralelo isomorfismo}; existence and
invertibility in the present case instead come from the linear ODE on the
Banach space proved in the preceding chapter.

If $\mathbf{T}_\gamma:=\dot\gamma$ denotes the tangent field, a curve $C^2$
is an \textbf{affine geodesic} if
$\nabla_t\mathbf{T}_\gamma=0$. In coordinates, this condition is
\begin{equation}\label{eq:geodesica-hilbert}
 \ddot\gamma+\Gamma_\gamma(\dot\gamma,\dot\gamma)=0.
\end{equation}
Every affine geodesic has constant speed, since
$\frac d{dt}|\mathbf{T}_\gamma|_{\mathbf{g}}^2
=2\mathbf{g}(\nabla_t\mathbf{T}_\gamma,\mathbf{T}_\gamma)=0$.
This definition and its local expression are the analogues of
Definition~\ref{def:geodesica}; the subsequent analytic novelty is solving the
equation in a Banach space, not the calculation that produces it.

\begin{definition}[Spray]
\label{def:spray-banach}
Let $M$ be a smooth Banach manifold.
\index{spray}
Let $h_t\colon TM\longrightarrow TM$ be the fiberwise homothety $h_t(v)=tv$. A
\textbf{spray} is a
smooth vector field $\mathbf{S}$ on $TM$ satisfying
\[
 T\pi\circ \mathbf{S}=\operatorname{id}_{TM},
 \qquad
 \mathbf{S}(h_t(v))=Th_t\bigl(t\mathbf{S}(v)\bigr)
\]
for every $t\in\mathbb R$ and $v\in TM$. In a chart it is written as
\[
 \mathbf{S}(x,v)=\bigl((x,v),(v,Q_x(v))\bigr),
 \qquad Q_x(tv)=t^2Q_x(v).
\]
Smoothness at the zero section implies that $Q_x(v)$ is the diagonal of a
continuous symmetric bilinear form. Indeed, for fixed $v$,
$Q_x(tv)=t^2Q_x(v)$ implies, by differentiating twice at $t=0$,
\[
 Q_x(v)=\frac12D_v^2Q_x(0)(v,v).
\]
The second derivative is symmetric, bilinear, and continuous; moreover, it depends
smoothly on $x$ in the operator norm by Fréchet calculus.
A field that satisfies only the first
condition is called a second-order field.
\end{definition}

The Levi--Civita spray is given by
\begin{equation}\label{eq:spray-levi-civita-hilbert}
 \mathbf{S}_{\mathbf{g}}(x,v)=\bigl((x,v),(v,-\Gamma_x(v,v))\bigr).
\end{equation}
Law \eqref{eq:ley-transformacion-christoffel-hilbert} is exactly what
ensures that these expressions glue together. A curve $\gamma$ satisfies
\eqref{eq:geodesica-hilbert} if and only if $t\mapsto\dot\gamma(t)$ is an integral
curve of $\mathbf{S}_{\mathbf{g}}$.

\begin{proposition}\label{prop:spray-exponencial-hilbert}
Let $(M,\mathbf{g})$ be a Hilbert manifold with a strong Riemannian metric.
For each $v\in T_xM$, there exists a unique maximal geodesic $\gamma_v$ with
$\gamma_v(0)=x$ and $\dot\gamma_v(0)=v$. Its domain is an open interval
containing $0$, and its dependence on $(t,v)$ is smooth on the maximal
domain. The set
\[
 \mathcal D:=\{v\in TM\mid\gamma_v\text{ is defined on }[0,1]\}
\]
is open, contains the zero section, and the map
\[
 \exp\colon \mathcal D\longrightarrow M,
 \qquad \exp(v):=\gamma_v(1)
\]
is smooth. Moreover, $d(\exp_x)_0=\operatorname{id}_{T_xM}$; therefore,
$\exp_x$ is a diffeomorphism between neighborhoods of the origin and $x$.
\end{proposition}

\begin{proof}
In coordinates, \eqref{eq:geodesica-hilbert} is equivalent to the first-order
system $(\dot x,\dot v)=(v,-\Gamma_x(v,v))$.
Theorem~\ref{teo:flujo-maximal-banach} applied to the spray gives existence,
uniqueness, openness of the maximal domain, and smooth dependence on the initial values.
Homogeneity of the spray and uniqueness give the scaling identity
\begin{equation}\label{eq:escala-geodesica-exponencial}
 \gamma_v(t)=\gamma_{tv}(1)=\exp_x(tv)
\end{equation}
whenever the two sides are defined. In particular,
\[
 d(\exp_x)_0(v)
 =\frac d{dt}\biggr|_{t=0}\exp_x(tv)
 =\dot\gamma_v(0)=v.
\]
The scaling equality is the same geometric step as in
Lemma~\ref{lema: reescalamiento geodesicas}.
Theorem~\ref{teo:funcion-inversa-local-banach} concludes that $\exp_x$ is a
diffeomorphism between neighborhoods of $0$ and $x$. In finite dimensions, the geodesic
field, maximal existence, and exponential map appear in
Proposition~\ref{prop: campo geodesico},
Theorem~\ref{teo: existencia unicidad geodesicas}, and
Definition~\ref{def:variedades-riemannianas-mapeo-exponencial}, with their
properties in Proposition~\ref{propiedades mapeo exponencial}. The differential
argument is identical; here we use the flow and inverse function theorem in Banach spaces.
\end{proof}

For an absolutely continuous curve whose velocity has integrable norm,
we define its length; if $|\dot\gamma|_{\mathbf{g}}^2$ is also integrable, we define its
energy by
\[
 L_{\mathbf{g}}(\gamma)=\int_a^b|\dot\gamma|_{\mathbf{g}}\,dt,
 \qquad
 E_{\mathbf{g}}(\gamma)=\frac12\int_a^b|\dot\gamma|_{\mathbf{g}}^2\,dt.
\]
The Cauchy--Schwarz inequality gives
\begin{equation}\label{eq:longitud-energia-hilbert}
 L_{\mathbf{g}}(\gamma)^2\leq2(b-a)E_{\mathbf{g}}(\gamma),
\end{equation}
with equality if and only if the speed is constant almost everywhere.

\begin{proposition}[First variation]\label{prop:primera-variacion-hilbert}
Let $(M,\mathbf{g})$ be a Hilbert manifold with a strong Riemannian metric and Levi--Civita connection.
Let $\alpha\colon (-\varepsilon,\varepsilon)\times[a,b]\longrightarrow M$ be a
smooth variation, $\mathbf{T}=\partial_t\alpha$, $\mathbf{V}=\partial_s\alpha$, and
$\gamma=\alpha_0$. On the following right-hand sides, $\mathbf{T}$ and $\mathbf{V}$ are evaluated
at $s=0$. Then
\[
 \frac d{ds}\biggr|_{s=0}E_{\mathbf{g}}(\alpha_s)
 =\bigl[\mathbf{g}(\mathbf{V},\mathbf{T})\bigr]_a^b-
   \int_a^b \mathbf{g}\bigl(\mathbf{V},\nabla_t\mathbf{T}\bigr)\,dt.
\]
The same formulas hold for variations of class $C^2$. If
$\mathbf{T}$ does not vanish along $\gamma$, we also have
\begin{equation}\label{eq:primera-variacion-longitud-hilbert}
 \frac d{ds}\biggr|_{s=0}L_{\mathbf{g}}(\alpha_s)
 =\left[\frac{\mathbf{g}(\mathbf{V},\mathbf{T})}{|\mathbf{T}|_{\mathbf{g}}}\right]_a^b
 -\int_a^b \mathbf{g}\left(\mathbf{V},\nabla_t\frac{\mathbf{T}}{|\mathbf{T}|_{\mathbf{g}}}\right)dt.
\end{equation}
Consequently, a curve $C^2$ is critical for energy with fixed
endpoints if and only if it is an affine geodesic. A regular curve is critical for
length with fixed endpoints if and only if its arc length
parametrization is an affine geodesic.
\end{proposition}

\begin{proof}
Fix $0<\delta<\varepsilon$ and restrict the variation to the compact
rectangle $[-\delta,\delta]\times[a,b]$. The function
$(s,t)\mapsto \mathbf{g}_{\alpha(s,t)}(\mathbf{T}(s,t),\mathbf{T}(s,t))$ is smooth
and, for a $C^2$ variation, is at least $C^1$. Differentiation under the
integral sign is therefore justified by uniform continuity of its
partial derivative. Metric compatibility gives
\[
 \frac d{ds}\biggr|_{s=0}E_{\mathbf{g}}(\alpha_s)
 =\int_a^b \mathbf{g}(\nabla_s\mathbf{T},\mathbf{T})\,dt.
\]
Absence of torsion and the equality
$[\boldsymbol{\partial}_s,\boldsymbol{\partial}_t]=0$ on the parameter domain imply
$\nabla_s\mathbf{T}=\nabla_t\mathbf{V}$. Using metric compatibility a second time,
\begin{align*}
 \int_a^b \mathbf{g}(\nabla_t\mathbf{V},\mathbf{T})\,dt
 =\int_a^b\left(\frac d{dt}\mathbf{g}(\mathbf{V},\mathbf{T})-\mathbf{g}(\mathbf{V},\nabla_t\mathbf{T})\right)dt=\bigl[\mathbf{g}(\mathbf{V},\mathbf{T})\bigr]_a^b
   -\int_a^b \mathbf{g}(\mathbf{V},\nabla_t\mathbf{T})\,dt,
\end{align*}
which is the energy formula, including its boundary terms.

If $\mathbf{T}$ does not vanish along $\gamma$, compactness of $[a,b]$ and continuity allow us to
shrink $\varepsilon$ so that $|\mathbf{T}|_{\mathbf{g}}$ remains bounded away from zero. We can
again differentiate under the integral sign, and
\[
 \partial_s|\mathbf{T}|_{\mathbf{g}}
 =\frac{\mathbf{g}(\nabla_s\mathbf{T},\mathbf{T})}{|\mathbf{T}|_{\mathbf{g}}}
 =\mathbf{g}\left(\nabla_t\mathbf{V},\frac{\mathbf{T}}{|\mathbf{T}|_{\mathbf{g}}}\right).
\]
The identity
\[
 \frac d{dt}\mathbf{g}\left(\mathbf{V},\frac{\mathbf{T}}{|\mathbf{T}|_{\mathbf{g}}}\right)
 =\mathbf{g}\left(\nabla_t\mathbf{V},\frac{\mathbf{T}}{|\mathbf{T}|_{\mathbf{g}}}\right)
  +\mathbf{g}\left(\mathbf{V},\nabla_t\frac{\mathbf{T}}{|\mathbf{T}|_{\mathbf{g}}}\right)
\]
and integration over $[a,b]$ prove
\eqref{eq:primera-variacion-longitud-hilbert}.

To justify the variational consequence without assuming that every field along
$\gamma$ comes from a global variation, fix a compact subinterval
$J\subset(a,b)$ whose image is contained in a chart
$\varphi\colon U\longrightarrow H$. If $\mathbf{W}$ is a smooth field along
$\gamma$, with support contained in the interior of $J$, for small $|s|$
we define on $J$
\[
 \alpha(s,t)
 :=\varphi^{-1}\!\left(
   \varphi(\gamma(t))+s\,d\varphi_{\gamma(t)}\mathbf{W}(t)\right)
\]
and set $\alpha(s,t)=\gamma(t)$ outside $J$. Compactness of the support
allows the same interval of values of $s$ to be chosen, and vanishing of $\mathbf{W}$ in
a neighborhood of the endpoints of $J$ makes the gluing smooth. If the central curve
is only $C^2$, take $C^2$ fields along it, and the same gluing
produces a $C^2$ variation; the first variation formula remains
applicable. This variation has
variational field $\mathbf{W}$ and fixed endpoints.

If the first variation of energy vanishes for all these variations,
fix $z\in H$ and take
\[
 \mathbf{W}(t)=\chi(t)(d\varphi_{\gamma(t)})^{-1}z,
 \qquad \chi\in C_c^\infty(\mathring J).
\]
The formula already proved and the fundamental lemma of the calculus of variations
imply
\[
 \mathbf{g}_{\gamma(t)}\bigl((d\varphi_{\gamma(t)})^{-1}z,\nabla_t\mathbf{T}(t)\bigr)=0,
 \qquad t\in\mathring J.
\]
Since $z$ is arbitrary and $d\varphi_{\gamma(t)}$ is an isomorphism,
nondegeneracy of $\mathbf{g}$ gives $\nabla_t\mathbf{T}=0$. For length, the same argument
applied to $\nabla_t\bigl(\frac{\mathbf{T}}{|\mathbf{T}|_{\mathbf{g}}}\bigr)$ gives
$\nabla_t\bigl(\frac{\mathbf{T}}{|\mathbf{T}|_{\mathbf{g}}}\bigr)=0$. Upon parametrizing by arc length, this last field
is the tangent vector, and the identity becomes the affine geodesic
equation. The converses follow directly from the variation formulas.

Metric compatibility and differentiation along curves are those of
Definition~\ref{def conexion compatible y libre de torsion} and
Theorem~\ref{teo: derivada covariante a lo largo de una curva}; the differentiation and
integration argument is the same in finite dimensions. The
construction of variations uses a Hilbert chart and cutoff functions
on the parameter interval; it therefore remains valid even when the manifold
is not locally compact.
\end{proof}

For $p,q$ in the same connected component, let
\[
 d_{\mathbf{g}}(p,q):=\inf\{L_{\mathbf{g}}(\gamma)\mid
 \gamma\text{ is a piecewise $C^1$ curve from $p$ to $q$}\}.
\]

\begin{proposition}\label{prop:distancia-riemann-fuerte-topologia}
On each connected component of a strong Riemannian Hilbert manifold,
$d_{\mathbf{g}}$ is a finite distance inducing the manifold topology.
\end{proposition}

\begin{proof}
Fix $x_0$ and a chart
$\varphi\colon U\longrightarrow H$ with $\varphi(x_0)=0$. By
Proposition~\ref{prop:caracterizaciones-metrica-fuerte}, after shrinking
$U$ there exist $c,C>0$ such that
\[
 c\|w\|_H\leq |d\varphi_x^{-1}w|_{\mathbf{g},x}\leq C\|w\|_H
\]
for every $x\in U$ and $w\in H$. Choose $r>0$ so that
$B_H(0,3r)\subseteq\varphi(U)$ and set
$V:=\varphi^{-1}(B_H(0,\frac{r}{2}))$. If $x,y\in V$, the segment joining their
coordinates remains in $B_H(0,\frac{r}{2})$, and its lift has length at
most $C\|\varphi(y)-\varphi(x)\|_H$.

Now let $\gamma$ be a piecewise $C^1$ curve from $x$ to $y$. As long as
$\gamma$ remains in $\varphi^{-1}(B_H(0,2r))$, the chain rule gives
\begin{equation}\label{eq:primera-salida-distancia-hilbert}
 c\|\varphi(\gamma(t_1))-\varphi(\gamma(t_0))\|_H
 \leq L_{\mathbf{g}}\bigl(\gamma\restriction_{[t_0,t_1]}\bigr).
\end{equation}
We justify that a curve leaving this ball reaches its coordinate
sphere before leaving the chart. Let $t_1$ be the endpoint of the initial
interval during which the curve remains in the ball. For $s<t<t_1$,
\eqref{eq:primera-salida-distancia-hilbert} bounds the difference of its
coordinates by $c^{-1}$ times the length of the segment $[s,t]$. That length
tends to zero as $s,t\to t_1^{-}$, because the speed is integrable.
The coordinates therefore have a limit in the complete space $H$, of
norm at most $2r$. This limit belongs to $\varphi(U)$, since
$\overline{B_H(0,2r)}\subset B_H(0,3r)\subset\varphi(U)$.
Continuity of $\varphi^{-1}$ and $\gamma$, together with the Hausdorff
property of $M$, identifies $\gamma(t_1)$ with the preimage of the limit.
If an exit occurs, its norm must be $2r$: a smaller norm
would allow the interval to be extended inside the ball.

If the curve does not leave this ball, taking its endpoints gives
$L_{\mathbf{g}}(\gamma)\geq c\|\varphi(y)-\varphi(x)\|_H$. If it does, continuity
of $\|\varphi\circ\gamma\|_H$ gives a first time $t_1$ at which
its value is $2r$. Since $\|\varphi(x)\|_H<\frac{r}{2}$, inequality
\eqref{eq:primera-salida-distancia-hilbert} implies
\[
 L_{\mathbf{g}}(\gamma)
 \geq c\bigl(2r-\|\varphi(x)\|_H\bigr)
 >\frac{3cr}{2}.
\]
On the other hand,
$\|\varphi(y)-\varphi(x)\|_H<r$, so this last bound is greater than
$c\|\varphi(y)-\varphi(x)\|_H$. Taking the infimum over all curves
gives
\[
 c\|\varphi(y)-\varphi(x)\|_H
 \leq d_{\mathbf{g}}(x,y)
 \leq C\|\varphi(y)-\varphi(x)\|_H
\]
for $x,y\in V$. To control also the points outside $V$,
apply the same exit argument to the ball of radius $r/2$, with initial
point $x_0$. Every curve from $x_0$ to a point of $M\setminus V$ has
length at least $cr/2$. Consequently,
$B_{d_{\mathbf g}}(x_0,cr/2)\subseteq V$.
The local upper inequality proves that convergence in coordinates
implies convergence in distance; this inclusion and the lower inequality
prove the converse. They also separate $x_0$ from every other
point. Symmetry follows by reversing curves and the triangle inequality
by concatenating them. Classes of points joined by piecewise $C^1$ curves are open; connectedness of the
component shows that there is only one and proves finiteness.

The corresponding finite-dimensional definition and metric conclusion are
found in
Definition~\ref{def:variedades-riemannianas-distancia-riemanniana-de-p-a-q}
and Theorem~\ref{variedad riemanniana espacio métrico}. The coordinate
segment gives the upper bound; the first exit argument gives the lower
bound through local comparison of norms, without requiring
compact balls.
\end{proof}

For a weak metric, in general one obtains only a pseudodistance:
positivity of the length of every curve does not prevent the existence of curves between
two distinct points with arbitrarily small lengths.

The local part of Riemannian geometry survives in the strong case; the
global part does not do so automatically.

\begin{proposition}[The surviving implication of Hopf--Rinow]
\label{prop:completo-metrico-implica-geodesico-hilbert}
Every metrically complete strong Riemannian manifold is geodesically
complete.
\end{proposition}

\begin{proof}
Let $\gamma\colon [0,b)\longrightarrow M$ be a maximal geodesic with
$b<\infty$. Its speed is constant, so
$d_{\mathbf{g}}(\gamma(s),\gamma(t))\leq |t-s||\dot\gamma(0)|_{\mathbf{g}}$. By completeness,
$\gamma(t)$ converges to a point $q$ as $t\to b^{-}$. Choose a chart
$\varphi\colon U\longrightarrow H$ about $q$. After shrinking it,
there exist $c,C,K>0$ such that
\[
 c\|w\|_H\leq|d\varphi_x^{-1}w|_{\mathbf{g},x}\leq C\|w\|_H,
 \qquad
 \|\Gamma_x\|_{\mathcal L_s^2(H;H)}\leq K
\]
for $x\in U$ and $w\in H$. For $t$ sufficiently close to $b$, the curve
remains in this chart. If its constant speed is $a$, then
$\|(\varphi\circ\gamma)'(t)\|_H\leq \frac{a}{c}$ and the geodesic equation gives
\[
 \|(\varphi\circ\gamma)''(t)\|_H\leq K(\frac{a}{c})^2.
\]
Thus, $(\varphi\circ\gamma)'(t)$ is Cauchy as $t\to b^{-}$ and converges
to some $v\in H$. In the coordinate trivialization, the state converges to
$(\varphi(q),v)$; therefore,
$(\gamma(t),\dot\gamma(t))$ converges to
$(q,(d\varphi_q)^{-1}v)\in TM$. Local existence for
the spray extends the geodesic beyond $b$, a contradiction.

In finite dimensions, this is one of the implications in
Theorem~\ref{teo:variedades-riemannianas-hopf-rinow}. Here we reuse only
constant speed and local existence for the ODE; we do not use compactness of
closed balls, which fails in an infinite-dimensional Hilbert space.
\end{proof}

\begin{example}[Completeness without a minimizing geodesic]
\label{ej:elipsoide-grossman}
Let $(e_n)_{n\in\mathbb N}$ be the standard basis of $\ell^2(\mathbb N)$, let $a_1=1$ and
$a_n=1+2^{-n}$ for $n\geq2$, and consider the ellipsoid
\[
 \mathcal E
 :=\left\{x=(x_n)\in\ell^2(\mathbb N)\ \middle|\
          \sum_{n=1}^{\infty}\frac{x_n^2}{a_n^2}=1\right\}
\]
with the metric induced by $\ell^2(\mathbb N)$. It is a Hilbert submanifold: the
function $\displaystyle f(x)=\displaystyle\sum_{n=1}^{\infty}x_n^2/a_n^2$ is a continuous quadratic
form, hence smooth, and $Df(x)[x]=2$ on the level $f=1$.
Its differential there is surjective onto $\mathbb R$ and has complemented
kernel. The regular value theorem in
Chapter~\ref{cap:variedades-banach} provides the submanifold
structure. The induced metric is strong, because it uses the norm of
$\ell^2(\mathbb N)$ on each closed tangent subspace. The operator
$A(x_n)=(a_nx_n)$ is a bounded linear isomorphism, maps the unit sphere
onto $\mathcal E$, and satisfies
\[
 \|v\|_{\ell^2(\mathbb N)}\leq\|Av\|_{\ell^2(\mathbb N)}
 \leq(1+2^{-2})\|v\|_{\ell^2(\mathbb N)}.
\]
If $d_{\mathbb S}$ is the intrinsic distance on the unit sphere and
$d_{\mathcal E}$ that on $\mathcal E$, comparison of lengths gives, for
$x,y$ on the sphere,
\[
 d_{\mathbb S}(x,y)
 \leq d_{\mathcal E}(Ax,Ay)
 \leq\frac54d_{\mathbb S}(x,y).
\]
The Hilbert sphere is complete: the ambient distance satisfies
$\|x-y\|_{\ell^2(\mathbb N)}\leq d_{\mathbb S}(x,y)$ and, if
$\|x-y\|_{\ell^2(\mathbb N)}<1$, the great-circle arc contained in
$\operatorname{span}\{x,y\}$ has length
$2\arcsin\!\left(\frac{\|x-y\|_{\ell^2(\mathbb N)}}{2}\right)\leq2\|x-y\|_{\ell^2(\mathbb N)}$. Hence an intrinsic
Cauchy sequence converges first in $\ell^2(\mathbb N)$ to a point of the sphere and
then in the intrinsic distance. The preceding inequalities prove that
$\mathcal E$ is metrically complete, and
Proposition~\ref{prop:completo-metrico-implica-geodesico-hilbert} proves its
geodesic completeness.

Every piecewise $C^1$ curve $c\colon[0,1]\to\mathbb S$ joining $e_1$
to $-e_1$ has length at least $\pi$. To see this, set
$a(t)=\langle c(t),e_1\rangle$. Since $\langle\dot c(t),c(t)\rangle=0$,
\[
 |a'(t)|=|\langle\dot c(t),e_1-a(t)c(t)\rangle|
 \leq\|\dot c(t)\|_{\ell^2(\mathbb N)}\sqrt{1-a(t)^2}.
\]
For $0<\delta<1$, the function
$\theta_\delta(t)=\arccos((1-\delta)a(t))$ satisfies
$|\theta_\delta'(t)|\leq\|\dot c(t)\|_{\ell^2(\mathbb N)}$ on each segment.
Integrating and letting $\delta\to 0^{+}$, we obtain $\pi\leq L(c)$.
Moreover,
\[
 \|A\dot c(t)\|_{\ell^2(\mathbb N)}^2-\|\dot c(t)\|_{\ell^2(\mathbb N)}^2
 =\sum_{n=2}^{\infty}(a_n^2-1)|\langle\dot c(t),e_n\rangle|^2\geq0.
\]
If $L(A\circ c)=L(c)$, the nonnegative difference of the speeds
vanishes almost everywhere. Since $a_n>1$ for each $n\geq2$, each
component $\langle\dot c(t),e_n\rangle$ with that index vanishes almost
everywhere; the fundamental theorem of calculus gives
$\langle c(t),e_n\rangle=0$ for every $t\in[0,1]$. The curve
would remain in $\mathbb S\cap\mathbb Re_1=\{e_1,-e_1\}$, where there is no
continuous curve between those endpoints. Hence, $L(A\circ c)>L(c)$. In contrast, the semicircles
\[
 c_n(t)=\cos(\pi t)e_1+\sin(\pi t)e_n
\]
are considered for $n\geq2$ and $t\in[0,1]$, and satisfy
\[
 L(A\circ c_n)
 =\pi\int_0^1\sqrt{\sin^2(\pi t)+a_n^2\cos^2(\pi t)}\,dt
 \longrightarrow\pi.
\]
The integrand converges uniformly, since $a_n\to1$; this also
follows by dominated convergence with bound $5/4$. Thus
$d_{\mathcal E}(e_1,-e_1)=\pi$, but no curve realizes that distance. This
is Grossman's ellipsoid \cite{Grossman1965}.
\end{example}

\begin{remark}[Limitations of Hopf--Rinow]
\label{obs:hopf-rinow-infinito}
Proposition~\ref{prop:completo-metrico-implica-geodesico-hilbert} does not
restore the finite-dimensional equivalences. Metric completeness does not
imply that every pair of points can be joined by a minimizing geodesic, as
Example~\ref{ej:elipsoide-grossman} shows. Nor does metric completeness
force $\exp_p$ to be surjective: this is Atkin's counterexample
in \cite{Atkin1975}. The converse implication from global geodesic
completeness to metric completeness also fails, but this belongs to
a later work of Atkin.\footnote{C.~J. Atkin, \emph{Geodesic and metric
completeness in infinite dimensions}, Hokkaido Mathematical Journal
\textbf{26} (1997), 1--61,
\href{https://doi.org/10.14492/hokmj/1351257804}{doi:10.14492/hokmj/1351257804}.}
These assertions distinguish three properties: continuation of every
geodesic, convergence of Cauchy sequences, and existence of
minimizing curves. In infinite dimensions one lacks the local compactness
of closed balls that relates them in the classical theorem.

There are important partial results. Ekeland proved that, for fixed $p$ on a
strong, connected, metrically complete manifold, the set of points that can
be joined to $p$ by a unique minimizing geodesic contains a
dense $G_\delta$ \cite[Thm.~B]{Ekeland1978}. Under global hypotheses of
nonpositive curvature, conclusions of Cartan--Hadamard type
can also be recovered. None of these assertions justifies using the
finite-dimensional Hopf--Rinow theorem without checking its additional hypotheses.
\end{remark}

\section{Banach--Finsler structures}
\label{sec:banach-finsler}

\begin{definition}\label{def:estructura-banach-finsler}\glsadd{banach-finsler}
\index{Banach--Finsler structure}\index{norm!Finsler}
\index{Finsler structure!reversible}
Let $M$ be a Banach manifold of class at least $C^1$, modeled on $X$. A
\textbf{Banach--Finsler structure} is a family of norms
$F_x\colon T_xM\longrightarrow[0,\infty)$ satisfying:
\begin{enumerate}[label=\textup{(\roman*)}]
 \item $F_x$ is a complete norm inducing the given topology of $T_xM$
       for every $x$;
 \item $(x,v)\mapsto F_x(v)$ is continuous on $TM$;
 \item for each $x_0$ there is a chart
       $\varphi\colon U\longrightarrow X$ such that, if
       $\|w\|_x:=F_x((d\varphi_x)^{-1}w)$, for every $\varepsilon>0$ we can
       shrink $U$ so that
 \begin{equation}\label{eq:uniformidad-local-finsler}
 (1+\varepsilon)^{-1}\|w\|_{x_0}\leq \|w\|_x
 \leq(1+\varepsilon)\|w\|_{x_0},
       \qquad x\in U,\ w\in X.
 \end{equation}
\end{enumerate}
The structure $F$ is called \textbf{reversible} if $F_x(-v)=F_x(v)$ for
every $x\in M$ and every $v\in T_xM$. This property concerns the structure
$F$, not the Banach manifold without that structure. Here ``norm'' has
its usual meaning, so absolute homogeneity implies this
equality: all Banach--Finsler structures in this definition are
reversible. Differentiability of $F$ at the zero section is not required.
Classical nonreversible Finsler metrics will be introduced later
through Minkowski functionals and are not instances of this definition.
\end{definition}

Condition \eqref{eq:uniformidad-local-finsler}, sometimes called
\textbf{admissibility}, is the Palais condition. It is stronger than requiring that,
for each $x$, there exist constants $c_x,C_x>0$ such that
$c_x\|w\|_X\leq F_x((d\varphi_x)^{-1}w)\leq C_x\|w\|_X$ for every $w\in X$:
in infinite dimensions, pointwise continuity in $w$ does not by
itself provide uniformity on the unit ball. Taking
$\|\cdot\|_X=\|\cdot\|_{x_0}$ and $C=1+\varepsilon$ gives
\begin{equation}\label{eq:comparacion-local-finsler-C}
 C^{-1}\|w\|_X\leq\|w\|_x\leq C\|w\|_X,
 \qquad x\in U,\quad w\in X.
\end{equation}

The definition does not depend on the chosen admissible chart. Indeed, if
$\psi$ is another chart and $\tau=\psi\circ\varphi^{-1}$, then
\[
 \|z\|_{x,\psi}
 =\bigl\|D\tau(\varphi(x))^{-1}z\bigr\|_{x,\varphi}.
\]
Set
\[
 A_x:=D\tau(\varphi(x))^{-1}D\tau(\varphi(x_0)).
\]
Continuity of inversion in the general linear group implies
$\|A_x-I\|_{\mathcal L(X)}\to0$, where the operator norm is calculated with
$\|\,\cdot\,\|_{x_0,\varphi}$. Given $\eta>0$, choose
$\varepsilon,\delta>0$ so that
\[
 (1+\varepsilon)(1+\delta)\leq1+\eta,
 \qquad
 (1+\varepsilon)^{-1}(1-\delta)\geq(1+\eta)^{-1}.
\]
After shrinking the domain, admissibility in the chart $\varphi$ holds with
$\varepsilon$ and $\|A_x-I\|\leq\delta$. For
$w=D\tau(\varphi(x_0))^{-1}z$ we then obtain
\begin{align*}
 \|z\|_{x,\psi}=\|A_xw\|_{x,\varphi}
 &\leq(1+\varepsilon)(1+\delta)\|w\|_{x_0,\varphi},\\
 \|z\|_{x,\psi}
 &\geq(1+\varepsilon)^{-1}(1-\delta)\|w\|_{x_0,\varphi}.
\end{align*}
Since $\|w\|_{x_0,\varphi}=\|z\|_{x_0,\psi}$, this is
\eqref{eq:uniformidad-local-finsler} in the chart $\psi$, with constant
$1+\eta$.

\begin{proposition}\label{prop:riemann-fuerte-es-banach-finsler}
If $\mathbf{g}$ is a strong Riemannian metric on a Hilbert manifold, then
$F_x(v):=|v|_{\mathbf{g},x}$ is a Banach--Finsler structure. Its length and
distance agree with $L_{\mathbf{g}}$ and $d_{\mathbf{g}}$.
\end{proposition}

\begin{proof}
Only admissibility needs to be checked. In a chart, fix $x_0$ and
conjugate $G(x)$ by $G(x_0)^{-\frac{1}{2}}$. The positive operator
\[
 B(x):=G(x_0)^{-\frac{1}{2}}G(x)G(x_0)^{-\frac{1}{2}}
\]
converges to $I$ in norm as $x\to x_0$. For $w\in H$, set
$z=G(x_0)^{\frac{1}{2}}w$. Then
$\|w\|_{x_0}^2=\|z\|_H^2$ and
$\|w\|_x^2=\langle B(x)z,z\rangle_H$. If
$\|B(x)-I\|<\delta<1$, it follows that
\[
 (1-\delta)\|w\|_{x_0}^2
 \leq\|w\|_x^2
 \leq(1+\delta)\|w\|_{x_0}^2.
\]
Given $\varepsilon>0$, a choice of $\delta>0$ satisfying
$1+\delta\leq(1+\varepsilon)^2$ and
$1-\delta\geq(1+\varepsilon)^{-2}$ gives
\eqref{eq:uniformidad-local-finsler}. Equality of lengths and
distances follows from the definitions.
\end{proof}

\begin{definition}\label{def:norma-dual-finsler}
Let $(M,F)$ be a Banach--Finsler manifold.
\index{dual norm!Finsler}
The \textbf{dual Finsler norm} of $\xi\in T'_xM$ is
\[
 F_x^*(\xi):=\sup\{ |\xi(v)|\mid v\in T_xM,\ F_x(v)\leq1\}.
\]
For $I\in C^1(M,\mathbb R)$ we write
$\|dI(x)\|_{F,*}:=F_x^*(dI(x))$.
\end{definition}

\begin{proposition}\label{prop:dual-finsler-coordenadas}
Let $(M,F)$ be a Banach--Finsler manifold modeled on $X$.
In a chart satisfying \eqref{eq:uniformidad-local-finsler}, for
$\lambda\in X'$ and $x$ in a smaller domain where
\eqref{eq:comparacion-local-finsler-C} holds,
\[
 C^{-1}\|\lambda\|_{X'}
 \leq F_x^*(\lambda\circ d\varphi_x)
 \leq C\|\lambda\|_{X'},
\]
where the norm of $X'$ is dual to $\|\cdot\|_X=\|\cdot\|_{x_0}$.
In particular, $x\mapsto\|dI(x)\|_{F,*}$ is continuous for $I\in C^1(M)$.
\end{proposition}

\begin{proof}
From \eqref{eq:comparacion-local-finsler-C} we obtain the inclusions
\[
 C^{-1}B_X\subseteq B_x\subseteq CB_X,
 \qquad B_x:=\{w\in X\mid\|w\|_x\leq1\}.
\]
Here $B_X:=\{w\in X\mid\|w\|_X\leq1\}$ is the closed unit ball.
Taking suprema of $|\lambda|$ over these sets gives the
dual estimate.

For continuity, fix $x_0$ and $\eta>0$. Admissibility allows the chart to be
shrunk so that the dual norms at $x$ and $x_0$ differ by a
factor between $(1+\eta)^{-1}$ and $1+\eta$. If
$\lambda_x:=dI(x)\circ(d\varphi_x)^{-1}$, then
$\lambda_x\to\lambda_{x_0}$ in $X'$ and
\begin{align*}
 \bigl|\|\lambda_x\|_{x,*}-\|\lambda_{x_0}\|_{x_0,*}\bigr|
 \leq{}&\|\lambda_x-\lambda_{x_0}\|_{x,*}\\
 &+\bigl|\|\lambda_{x_0}\|_{x,*}
          -\|\lambda_{x_0}\|_{x_0,*}\bigr|.
\end{align*}
The precise dual comparison is
\[
 (1+\eta)^{-1}\|\mu\|_{x_0,*}
 \leq\|\mu\|_{x,*}
 \leq(1+\eta)\|\mu\|_{x_0,*},
 \qquad \mu\in X'.
\]
Therefore, the first term on the preceding right-hand side does not exceed
$(1+\eta)\|\lambda_x-\lambda_{x_0}\|_{x_0,*}$ and the second does not exceed
$\eta\|\lambda_{x_0}\|_{x_0,*}$. First choose $\eta$ and then $x$; both
terms can be made arbitrarily small. This proves continuity.
\end{proof}

For a piecewise $C^1$ curve, set
\[
 L_F(\gamma)=\int_a^bF_{\gamma(t)}(\dot\gamma(t))\,dt,
 \qquad
 d_F(x,y)=\inf\{L_F(\gamma)\mid
 \gamma\text{ is a piecewise $C^1$ curve from $x$ to $y$}\}.
\]
Length is invariant under increasing piecewise $C^1$
reparametrizations. Symmetry of the norms implies
$L_F(\gamma^{-1})=L_F(\gamma)$.

\begin{theorem}\label{teo:distancia-banach-finsler-topologia}
On each connected component of a Banach--Finsler manifold, $d_F$ is a
finite distance inducing the manifold topology. More precisely,
for each $x_0$ there exist a chart $\varphi\colon U\longrightarrow X$, a neighborhood
$V\subseteq U$ of $x_0$, and $C\geq1$ such that, for $x,y\in V$,
\[
 C^{-1}\|\varphi(y)-\varphi(x)\|_X
 \leq d_F(x,y)
 \leq C\|\varphi(y)-\varphi(x)\|_X.
\]
\end{theorem}

\begin{proof}
Translating coordinates, suppose that $\varphi(x_0)=0$. First shrink
$U$ so that \eqref{eq:comparacion-local-finsler-C} holds there, and then choose
$r>0$ so that $B_X(0,3r)\subseteq\varphi(U)$. Let
\[
 V:=\varphi^{-1}(B_X(0,\frac{r}{2})).
\]
If $x,y\in V$, the segment
$z(t)=(1-t)\varphi(x)+t\varphi(y)$ remains in $B_X(0,\frac{r}{2})$, and its lift
$\sigma=\varphi^{-1}\circ z$ satisfies
\[
 L_F(\sigma)
 \leq C\int_0^1\|z'(t)\|_X\,dt
 =C\|\varphi(y)-\varphi(x)\|_X.
\]
This is the upper bound.

Now let $\gamma$ be any curve from $x$ to $y$. While its image remains
in $\varphi^{-1}(B_X(0,2r))$, the chain rule and
\eqref{eq:comparacion-local-finsler-C} give
\begin{equation}\label{eq:cota-inferior-curva-carta-finsler}
 \|\varphi(\gamma(t_1))-\varphi(\gamma(t_0))\|_X
 \leq C L_F(\gamma\restriction_{[t_0,t_1]}).
\end{equation}
If $\gamma$ does not leave that ball, apply the inequality to its endpoints and
obtain
$L_F(\gamma)\geq C^{-1}\|\varphi(y)-\varphi(x)\|_X$. If it does, there is
a first time $t_1$ at which the norm of its coordinate equals $2r$.
Indeed, estimate
\eqref{eq:cota-inferior-curva-carta-finsler} makes the coordinates a
Cauchy family as the endpoint of the initial segment is approached.
Completeness of $X$ provides a limit in $\overline{B_X(0,2r)}$, which is
contained in $\varphi(U)$. Continuity of the curve and the Hausdorff
property of $M$ identify that limit with the coordinate of $\gamma(t_1)$.
If its norm were less than $2r$, the curve would remain in the ball a little
longer. This is the exit argument of
Proposition~\ref{prop:distancia-riemann-fuerte-topologia}, based on finite
length and completeness of the model. Since
$\|\varphi(x)\|_X<\frac{r}{2}$, the same inequality gives
\[
 L_F(\gamma)\geq C^{-1}\frac{3r}{2}
 \geq C^{-1}\|\varphi(y)-\varphi(x)\|_X,
\]
because the distance between two points of $B_X(0,\frac{r}{2})$ is less than $r$. Taking
the infimum over $\gamma$ gives the lower bound, even against curves
that leave the chart.

Nonnegativity, symmetry, and the triangle inequality
follow, respectively, from the definition, reversing curves, and
concatenating them. To include points outside $V$, the same argument
in the ball of radius $r/2$ proves that every curve from $x_0$ to
$M\setminus V$ has length at least $r/(2C)$. Thus,
$B_{d_F}(x_0,r/(2C))\subseteq V$. This inclusion and the lower bound inside
$V$ separate $x_0$ from all other points and show that a sequence
converging to $x_0$ in distance eventually enters $V$ and converges in
coordinates. The upper bound proves the converse. Finally, the relation
of accessibility by piecewise $C^1$ curves partitions $M$ into open
sets: from an accessible point one can reach every point of a sufficiently small
convex chart. By connectedness, each connected component is a
single accessibility class, so $d_F$ is finite on it. The two local
bounds show that the balls of $d_F$ and the coordinate neighborhoods generate the
same topology.
\end{proof}

\begin{definition}\label{def:completitud-finsler}
A Banach--Finsler manifold is \textbf{metrically complete} if each
component, equipped with the restriction of $d_F$, is a complete metric
space.
\end{definition}

\begin{definition}[$C^{1,1}$ regularity and regularity of the norm]
\label{def:finsler-C11}
Let $(M,F)$ be a Banach--Finsler manifold. We say it is of class
$C^{1,1}$ when $M$ is equipped with a $C^{1,1}$ atlas in the sense of
Definition~\ref{def:variedad-banach-C11}; regularity of coordinate
representations will be checked in this atlas. This condition does not
impose Lipschitz regularity on the family of norms $F_x$. A
Banach--Finsler structure is \textbf{locally Lipschitz in the base} if,
for each $x_0$, there exist a smaller admissible chart and $L>0$ such that
\begin{equation}\label{eq:norma-finsler-lipschitz-base}
 |\|w\|_x-\|w\|_y|
 \leq L\|\varphi(x)-\varphi(y)\|_X\,\|w\|_{x_0},
 \qquad x,y\in U,\quad w\in X.
\end{equation}
These two conditions are independent: the first concerns the atlas and the
second the family of norms. The second is expressed in charts of the chosen
$C^{1,1}$ atlas, since transformation of a tangent norm contains the
differential of the chart change.
\end{definition}

The term ``Finsler manifold of class $C^{1,1}$'' sometimes groups together
different hypotheses. The Palais--Struwe arguments use the
following ingredients:
\begin{enumerate}[label=\textup{(\alph*)}]
 \item $M$ is a paracompact Banach manifold with a $C^{1,1}$ atlas;
 \item $F$ is a continuous admissible structure in the sense of
       Definition~\ref{def:estructura-banach-finsler};
 \item the covers involved admit partitions of unity
       locally Lipschitz in the charts of the $C^{1,1}$ atlas.
\end{enumerate}
In this construction, the partitions in the last item do not constitute
an additional hypothesis on the model space. On each component,
Theorem~\ref{teo:distancia-banach-finsler-topologia} and
Theorem~\ref{teo:paracompacidad-espacios-metricos} allow
Lemma~\ref{lem:encogimiento-localmente-finito-banach} to be applied. Distances to
the complements of the shrunk open sets, divided by their positive,
locally finite sum, produce these partitions. The construction
is developed in Theorem~\ref{teo:existencia-pseudogradiente-finsler}.
When differentiable partitions are required, one instead uses
Theorem~\ref{teo:particiones-unidad-variedad-banach}, with the regularity and
model space hypotheses stated there.

Continuity and admissibility of $F$ suffice for the topology of $d_F$ and for
transporting strict inequalities from one fiber to nearby fibers. The
$C^{1,1}$ atlas and locally Lipschitz partitions are the ingredients that
allow local fields to be glued and an ODE with uniqueness to be obtained. The
additional condition \eqref{eq:norma-finsler-lipschitz-base} is invoked only when
the norm itself needs to vary with Lipschitz control. None of these
hypotheses asserts that $v\mapsto F_x(v)^2$ has a vertical Hessian or that a
Chern connection exists. This is the setting of the deformation theory of Palais
and Struwe \cite[\S II.3.7]{Struwe2008}.

\begin{example}[Complete submanifolds of a Banach space]
\label{ej:subvariedad-banach-finsler-completa}
Let $N$ be a closed split submanifold of a Banach space $X$. The norm
induced by $X$ on each $T_xN$ defines a complete Banach--Finsler structure.
In a local parametrization
$\theta\colon V\subseteq X_0\longrightarrow N$, continuity of
$D\theta$ in the operator norm and the fact that $D\theta(y_0)$ is an
isomorphism onto its image allow us to use
$\|w\|_0:=\|D\theta(y_0)w\|_X$. This norm is complete because
$D\theta(y_0)$ is a topological isomorphism onto the closed subspace
$T_{\theta(y_0)}N$. Given $\varepsilon>0$, shrink $V$ so that
\[
 \|D\theta(y)-D\theta(y_0)\|_{0,X}
 <\frac{\varepsilon}{1+\varepsilon}.
\]
Here $\|\,\cdot\,\|_{0,X}$ is the operator norm from the space
$(X_0,\|\,\cdot\,\|_0)$ to $(X,\|\,\cdot\,\|_X)$.
The triangle inequality and its reverse form then give
\[
 (1+\varepsilon)^{-1}\|w\|_0
 \leq\|D\theta(y)w\|_X
 \leq(1+\varepsilon)\|w\|_0,
\]
which is exactly admissibility of the induced norms. To prove
completeness, first observe
that, for every curve $\gamma$ in $N$,
\[
 \|\gamma(b)-\gamma(a)\|_X\leq L_F(\gamma),
\]
so an intrinsic Cauchy sequence is Cauchy in $X$. Its limit
belongs to $N$ because $N$ is closed. In an adapted submanifold chart,
the local parametrization and its differential are bounded after shrinking
the chart. The coordinate segment then gives
$d_F(x,p)\leq C_0\|\psi(x)-\psi(p)\|_{X_0}$ for $x$ close to the limit $p$,
where $C_0>0$ depends on the smaller chart. Therefore, convergence in $X$,
and with it convergence in the chart, implies
intrinsic convergence. Global closedness and local control of lengths
play different roles in this argument.
\end{example}

\begin{proposition}[Continuation criterion]\label{prop:flujo-acotado-finsler}
Let $(M,F)$ be a complete Banach--Finsler manifold, and let $\mathbf{X}$ be a field defined
on all of $M$, continuous and locally Lipschitz in coordinates. Suppose the
atlas of $M$ is at least $C^{1,1}$. If there exists $A\geq0$ such that
$F_x(\mathbf{X}(x))\leq A$ for every $x$, every maximal integral curve of $\mathbf{X}$ is
defined on $\mathbb R$. The global flow is continuous, and its time maps are
homeomorphisms with inverse $\Phi_{-t}$. If the atlas is $C^{r+1}$ and
$\mathbf{X}$ is $C^r$, with $r\in\mathbb N$, the flow and its time maps
are $C^r$; if the atlas and field are smooth, so is the flow.
\end{proposition}

\begin{proof}
If the maximal right endpoint were $b<\infty$, then, for $s<t<b$,
$d_F(\gamma(s),\gamma(t))
\leq L_F(\gamma\restriction_{[s,t]})\leq A(t-s)$. The curve is
Cauchy and, by completeness, converges to a point $p$. Since $d_F$ induces the
manifold topology, the tail of $\gamma$ enters a chart about
$p$. Remark~\ref{obs:campos-flujos-atlas-C11}, applied with
initial value $p$ at time $b$, produces a local solution. To check gluing at $b$, write the integral equation in the
final chart and let the upper endpoint tend to $b$. The field is
continuous and the curve has a limit, so the integral equation also holds
with the limiting initial value; its derivative at $b$ is the field value at $p$.
Local uniqueness identifies this continuation with the solution starting at
$p$. This gives an integral curve beyond $b$,
contradicting maximality. The left endpoint is treated in the same way, or
by applying the argument to $-\mathbf{X}$.

The chart argument of Theorem~\ref{teo:flujo-maximal-banach} remains
applicable: changes in a $C^{1,1}$ atlas preserve the local
Lipschitz condition on fields, by the estimate in
Section~\ref{sec:campos-curvas-integrales-banach}. The continuous dependence
theorem in the appendix gives continuity of the flow. If the atlas is $C^{r+1}$,
its tangent changes are $C^r$, and the same construction gives the
asserted $C^r$ dependence. Uniqueness gives
$\Phi_{t+s}=\Phi_t\circ\Phi_s$ and $\Phi_t^{-1}=\Phi_{-t}$; $C^r$ dependence
on the initial values gives the last assertion.
\end{proof}

\begin{remark}[Gradient versus pseudogradient]
For a strong metric, $\operatorname{grad}I=\sharp_{\mathbf{g}}(dI)$ and
$dI(\operatorname{grad}I)=|\operatorname{grad}I|_{\mathbf{g}}^2$. In a Banach--Finsler
structure there is generally no canonical choice of a vector that
realizes the dual norm, and in a nonreflexive Banach space the supremum may not
be attained. At a point where $dI\neq0$, deformation only requires choosing
locally, by the definition of the
supremum, directions $v$ with $F(v)\leq1$ and
$dI(v)>\theta\|dI\|_{F,*}$, $0<\theta<1$. One may also leave a
strict margin in the norm. Fixing $x_0$, choose
$\theta<\theta_1<1$ and $v_0\in T_{x_0}M$ with $F_{x_0}(v_0)\leq1$ and
$dI_{x_0}(v_0)>\theta_1\|dI(x_0)\|_{F,*}$. Multiply $v_0$ by
$a\in(\theta/\theta_1,1)$. Then $F_{x_0}(av_0)<1$ and
$dI_{x_0}(av_0)>\theta\|dI(x_0)\|_{F,*}$.
Extend $av_0$ as a field with constant components in a chart.
Continuity of $F$, $dI$, and the dual norm preserves both
strict inequalities in a neighborhood of $x_0$. In a $C^{1,1}$ atlas,
this field is locally Lipschitz in the other charts as well, without
having to divide it by a norm that is only continuous. To glue the directions without losing this inequality, one needs a
nonnegative locally finite partition of unity of the required
regularity. To obtain a locally Lipschitz pseudogradient, one uses
the partitions constructed from distances in
Theorem~\ref{teo:existencia-pseudogradiente-finsler}; their existence rests
on Theorem~\ref{teo:paracompacidad-espacios-metricos} and
Lemma~\ref{lem:encogimiento-localmente-finito-banach}. The model space need not
admit differentiable bump functions.
Theorem~\ref{teo:particiones-unidad-variedad-banach} concerns the
different problem of obtaining differentiable partitions. None of
these constructions enters the definition of $F$ or the
proof that $d_F$ induces the topology. This is the functional-analytic reason for
the pseudogradient vector fields used in variational deformations.
\end{remark}

\section{Classical Finsler geometry in finite dimensions}
\label{sec:finsler-clasico}

\subsection{Minkowski functionals and vertical tensors}

\begin{definition}[Minkowski functional]
\label{def:funcional-minkowski}
\index{Minkowski functional}
Let $V$ be a finite-dimensional real vector space. A
\textbf{Minkowski functional} is a continuous function
$F\colon V\longrightarrow[0,\infty)$,
smooth on $V\setminus\{0\}$, positive away from $0$, and positively homogeneous,
$F(\lambda v)=\lambda F(v)$ for $\lambda>0$, such that
\[
 \mathbf{g}_v(u,w):=\frac12D^2(F^2)_v(u,w)
\]
is an inner product for every $v\neq0$. It is \textbf{reversible} if
$F(-v)=F(v)$.
\end{definition}

\begin{definition}\label{def:finsler-clasico}\glsadd{finsler-clasico}
\index{Finsler metric!classical}\index{fundamental tensor!Finsler}
On a finite-dimensional smooth manifold $M$, a \textbf{classical Finsler
metric} is a function $F\colon TM\longrightarrow[0,\infty)$, continuous at the
zero section and
smooth on the slit bundle
\[
 0_M:=\{0_x\in T_xM\mid x\in M\},
 \qquad
 \widetilde{TM}:=TM\setminus0_M,
\]
such that each restriction $F_x:=F\restriction_{T_xM}$ is a Minkowski
functional. In coordinates $(x^1,\ldots,x^n,v^1,\ldots,v^n)$, where
$n=\dim M$, the indices $i,j,k$ range over $\{1,\ldots,n\}$. The
\textbf{fundamental tensor} and \textbf{Cartan tensor} are
\[
 g_{ij}(x,v)=\frac12\frac{\partial^2 F^2}{\partial v^i\partial v^j}(x,v),
 \qquad
 C_{ijk}(x,v)=\frac14
 \frac{\partial^3F^2}{\partial v^i\partial v^j\partial v^k}(x,v).
\]
Intrinsically,
$\mathbf{g}_v=\frac12D_v^2(F_x^2)$ and
$\mathbf{C}_v=\frac14D_v^3(F_x^2)$, where the derivatives
$D_v$ and $d(F_x)_v$ appearing in this section are always vertical; that
is, they are computed in the fixed fiber $T_xM$.
\end{definition}

\begin{proposition}[Euler identities and the Legendre transformation]
\label{prop:euler-legendre-finsler}
Let $(M,F)$ be a finite-dimensional classical Finsler manifold.
For $v\neq0$ and $u,w\in T_xM$ we have
\begin{align*}
 d(F_x)_v(u)&=\frac{\mathbf{g}_v(v,u)}{F(v)},
 &\mathbf{g}_v(v,v)&=F(v)^2,
 &\mathbf{C}_v(v,u,w)&=0.
\end{align*}
The Legendre transformation
\[
 \mathcal L_x\colon T_xM\setminus\{0\}\longrightarrow
 T_x^*M\setminus\{0\},
 \qquad
 \mathcal L_x(v):=d\left(\frac12F_x^2\right)_v
 =\mathbf{g}_v(v,\,\cdot\,),
\]
is a positively homogeneous diffeomorphism.
\end{proposition}

\begin{proof}
The function $F_x^2$ is homogeneous of degree two. Euler's theorem gives
$d(F_x^2)_v(v)=2F_x(v)^2$. Differentiating this identity in the direction $u$
gives
\[
 D^2(F_x^2)_v(v,u)=d(F_x^2)_v(u),
\]
which proves the first two formulas. Since
$\mathbf{g}_{\lambda v}=\mathbf{g}_v$, its radial derivative is zero, yielding
$\mathbf{C}_v(v,u,w)=0$.

The vertical derivative of $\mathcal L_x$ at $v$ is the isomorphism
$u\mapsto \mathbf{g}_v(u,\,\cdot\,)$, so $\mathcal L_x$ is a local
diffeomorphism. To justify injectivity also on segments passing through
the origin, fix a Euclidean norm $|\,\cdot\,|$ on $T_xM$ and set
$L=\frac12F_x^2$. Homogeneity and compactness of the Euclidean sphere
give $L(v)\leq C|v|^2$ and $\|dL_v\|\leq C'|v|$ for constants $C,C'>0$.
Hence $L$ is $C^1$ on all of $T_xM$, with $dL_0=0$.
If $u\ne v$, the function
$h(t):=L(u+t(v-u))$ has a continuous derivative on $[0,1]$ and
\[
 h''(t)=\mathbf g_{u+t(v-u)}(v-u,v-u)>0
\]
except, at most, at the unique parameter where the segment passes through zero.
On each side of that parameter, $h'$ is strictly increasing, and
continuity prevents a jump when crossing it. Therefore,
$h'(1)>h'(0)$, that is,
$\bigl(dL_v-dL_u\bigr)(v-u)>0$. This proves strict convexity of $L$
and injectivity of its differential. Given
$\xi\in T_x^*M\setminus\{0\}$, the function $\xi$ attains a positive maximum on
the compact sphere $F_x(v)=1$. This constraint is regular, since Euler's
identity gives $d(F_x)_v(v)=F_x(v)=1$.
Corollary~\ref{teo:multiplicadores-lagrange-banach} provides
$\xi=\lambda\,d(F_x)_v=\lambda\mathcal L_x(v)$ with $\lambda>0$; the last
equality uses $F_x(v)=1$, and positivity follows by evaluation at $v$. By
homogeneity, $\xi=\mathcal L_x(\lambda v)$; thus, $\mathcal L_x$ is
surjective. A bijection that is a local diffeomorphism has a smooth inverse.
\end{proof}

Strong convexity includes the triangle inequality. Indeed, the preceding
Euler identities and the quotient rule give, for $v\neq0$,
\[
 D^2F_v(w,w)
 =\frac{\mathbf{g}_v(w,w)F(v)^2-\mathbf{g}_v(v,w)^2}{F(v)^3}\geq0
\]
by Cauchy--Schwarz for $\mathbf{g}_v$. The restriction of $F$ to any line is
convex away from the possible parameters where it passes through the origin;
homogeneity describes the only missing point: on a line $tw$
through the origin, the left slope is $-F(-w)$ and the right slope
is $F(w)$. The former is smaller than the latter, so convexity
is preserved across the origin as well. Therefore, $F$
is convex on all of $V$, and convexity and homogeneity imply
$F(u+w)\leq F(u)+F(w)$. Thus, a reversible Minkowski functional is a
norm, whereas a nonreversible one is a directed norm.

The Cartan tensor measures directional dependence of the fundamental tensor.
If $\mathbf{g}_v$ does not depend on $v$, then $\mathbf{C}_v=0$. Conversely, if
$\dim M\geq2$, each slit fiber $T_xM\setminus\{0\}$ is connected and
$\mathbf{C}=0$ implies that $\mathbf{g}_v$ is independent of $v$; therefore, $F_x^2$ is a
quadratic form. In dimension one, the same conclusion additionally requires
reversibility: without it, the two rays of a slit fiber may carry
different quadratic forms. If $F$ is reversible, it satisfies
Definition~\ref{def:estructura-banach-finsler}. In
finite dimensions, fix a chart $\varphi$ and the norm
$\|w\|_0:=F_{x_0}((d\varphi_{x_0})^{-1}w)$. On the compact sphere
$\{w\in\mathbb R^n\mid\|w\|_0=1\}$, uniform continuity of
$F_x((d\varphi_x)^{-1}w)$ implies that, given $\varepsilon>0$, after
shrinking the chart,
\[
 (1+\varepsilon)^{-1}\|w\|_0
 \leq F_x((d\varphi_x)^{-1}w)
 \leq(1+\varepsilon)\|w\|_0
\]
for every $w$ and every $x$ in the smaller chart. This is admissibility;
the compactness used in the argument is exactly what is missing from the unit
ball of an infinite-dimensional Banach model.

\begin{example}[Minkowski spaces]
\label{ej:espacio-minkowski-finsler}
A finite-dimensional vector space equipped with a constant Minkowski
functional is a Finsler manifold. Its geodesics are straight lines, and the
nonlinear and Chern connections introduced below have
zero coefficients in linear coordinates. The Cartan tensor may be
nonzero: vanishing of the connection expresses that the geometry is locally
Minkowski, not that the functional is Euclidean.
\end{example}

\subsection{Spray, nonlinear connection, and Chern connection}

Write $n:=\dim M\geq1$. All Latin and Greek coordinate indices in
this section range over $\{1,\ldots,n\}$ unless another set is
specified. The matrix $(g^{ij})_{i,j=1}^n$ is the inverse of
$(g_{ij})_{i,j=1}^n$. Free indices remain fixed in each equality;
summed indices will be indicated in the corresponding sum.

For $L=\tfrac12F^2$, the extremals satisfy
\begin{equation}\label{eq:euler-lagrange-finsler}
 \frac d{dt}\frac{\partial L}{\partial v^i}(\gamma,\dot\gamma)
 -\frac{\partial L}{\partial x^i}(\gamma,\dot\gamma)=0.
\end{equation}
Since $\frac{\partial^2L}{\partial v^i\partial v^j}=g_{ij}$, the preceding equation
becomes
\begin{equation}\label{eq:geodesica-spray-finsler}
 \ddot x^i+2G^i(x,\dot x)=0,
\end{equation}
where
\begin{equation}\label{eq:coeficientes-spray-finsler}
 G^i(x,v)
 =\frac14\sum_{\ell=1}^n g^{i\ell}(x,v)
 \left(
 \sum_{k=1}^n\frac{\partial^2F^2}{\partial x^k\partial v^\ell}(x,v)v^k
 -\frac{\partial F^2}{\partial x^\ell}(x,v)
 \right).
\end{equation}
The coefficients $G^i$ are smooth on $\widetilde{TM}$ and homogeneous of degree
two. Under a coordinate change $y=\tau(x)$ they satisfy
\begin{equation}\label{eq:transformacion-spray-finsler}
 \widetilde G^\alpha(y,D\tau(x)v)
 =\sum_{i=1}^n\frac{\partial\tau^\alpha}{\partial x^i}(x)G^i(x,v)
 -\frac12\sum_{i,j=1}^n\frac{\partial^2\tau^\alpha}
 {\partial x^i\partial x^j}(x)v^iv^j,
\end{equation}
as can be checked by transforming \eqref{eq:geodesica-spray-finsler}. Hence
they define the \textbf{geodesic spray on the slit tangent bundle}
\[
 \mathbf{S}_F(x,v)=\sum_{i=1}^n v^i\boldsymbol{\partial}_{i}
 -2\sum_{i=1}^n G^i(x,v)\boldsymbol{\partial}_{v^i}
\]
on $\widetilde{TM}$. It is a second-order field and satisfies the homogeneity
identity of Definition~\ref{def:spray-banach} for $\lambda>0$:
\[
 \mathbf{S}_F(h_\lambda(v))=Th_\lambda\bigl(\lambda \mathbf{S}_F(v)\bigr),
 \qquad v\in\widetilde{TM},\quad\lambda>0.
\]
The restriction to positive scalars is essential when $F$ is not reversible.
In general, $\mathbf{S}_F$ does not extend smoothly to the zero section; hence it is not a
smooth spray on all of $TM$ in the sense of that definition, and constant
curves are included separately.

\begin{definition}[Canonical nonlinear connection]
\label{def:conexion-no-lineal-finsler}
Let $(M,F)$ be a finite-dimensional classical Finsler manifold.
\index{nonlinear connection!Finsler}
The coefficients
\[
 N^i_j(x,v):=\frac{\partial G^i}{\partial v^j}(x,v)
\]
define the \textbf{canonical nonlinear connection}. Denote the restricted tangent projection by $\widetilde\pi\colon\widetilde{TM}\to M$.
The vertical distribution is
\[
 \mathcal V_{(x,v)}:=\ker T_{(x,v)}\widetilde\pi
 =\operatorname{span}\left\{\boldsymbol{\partial}_{v^i}\mid 1\leq i\leq n\right\},
\]
and the horizontal distribution is spanned by the fields
\[
 \boldsymbol{\delta}_{x^j}
 :=\boldsymbol{\partial}_{j}
   -\sum_{i=1}^n N^i_j\boldsymbol{\partial}_{v^i},
\]
and the adapted vertical coforms are written as
\[
 \boldsymbol{\delta v}^{\,i}:=\mathbf{d}v^i+\sum_{j=1}^n N^i_j\mathbf{d}x^j.
\]
In particular,
$\boldsymbol{\delta v}^{\,i}(\boldsymbol{\delta}_{x^j})=0$ and
$\boldsymbol{\delta v}^{\,i}(\boldsymbol{\partial}_{v^j})=\delta^i_j$. If
$\mathcal H_{(x,v)}$ denotes the subspace spanned by
$\boldsymbol{\delta}_{x^j}$, then
$T_{(x,v)}\widetilde{TM}=\mathcal H_{(x,v)}\oplus\mathcal V_{(x,v)}$.
\end{definition}

Law \eqref{eq:transformacion-spray-finsler}, differentiated with respect to $v$, gives
the explicit transformation
\begin{equation}\label{eq:transformacion-conexion-no-lineal-finsler}
 \sum_{\beta=1}^n\widetilde N^\alpha_\beta(y,D\tau(x)v)
 \frac{\partial y^\beta}{\partial x^j}
 =\sum_{i=1}^n\frac{\partial y^\alpha}{\partial x^i}N^i_j(x,v)
  -\sum_{k=1}^n\frac{\partial^2y^\alpha}{\partial x^j\partial x^k}v^k.
\end{equation}
This identity shows that the horizontal distribution is intrinsic, although
the coefficients $N^i_j$ are not those of a linear connection on $TM$. Homogeneity gives
$\displaystyle\sum_{j=1}^n N^i_jv^j=2G^i$.

The local flow of $\mathbf{S}_F$ produces a maximal geodesic $\gamma_v$ for each
$v\neq0$. Homogeneity gives
$\gamma_v(t)=\gamma_{tv}(1)$ for $t>0$ whenever both sides are defined.
Adding the constant geodesics gives a star-shaped domain
$\mathcal D_p\subseteq T_pM$ and the exponential map
\[
 \exp_p\colon \mathcal D_p\longrightarrow M,
 \qquad \exp_p(v):=\gamma_v(1).
\]
It is smooth on $\mathcal D_p\setminus\{0\}$, of class $C^1$ on $0$, and
\[
 d(\exp_p)_0=\operatorname{id}_{T_pM}.
\]
To quantify this assertion, choose coordinates
$\varphi\colon U\longrightarrow\mathbb R^n$ with $\varphi(p)=0$ and use
$d\varphi_p$ to identify $T_pM$ with $\mathbb R^n$. Let
$\widehat\exp_p:=\varphi\circ\exp_p\circ(d\varphi_p)^{-1}$. Take
$\rho>0$ with $\overline{B(0,\rho)}\subseteq\varphi(U)$. Compactness of
$\overline{B(0,\rho)}\times S^{n-1}$ provides numbers
$M_0,M_x,M_v\geq0$ bounding, respectively,
\[
 |G(x,u)|,
 \qquad \|D_xG(x,u)\|,
 \qquad \|D_vG(x,u)\|,
 \qquad |u|=1.
\]
By homogeneity, for $w\neq0$ these bounds become
\[
 |G(x,w)|\leq M_0|w|^2,
 \quad \|D_xG(x,w)\|\leq M_x|w|^2,
 \quad \|D_vG(x,w)\|\leq M_v|w|.
\]
Choose $0<r\leq1$ so that $2r<\rho$ and $8M_0r\leq\frac{1}{2}$. If
$0<|v|<r$, the solution $x_v$ of
$x_v''=-2G(x_v,x_v')$, $x_v(0)=0$, $x_v'(0)=v$ satisfies
$|x_v'(t)|\leq2|v|$ for $0\leq t\leq1$. Indeed, while this
bound holds,
\[
 |x_v'(t)-v|
 \leq2M_0\int_0^t|x_v'(s)|^2ds
 \leq8M_0|v|^2\leq\frac12|v|,
\]
which improves the initial bound and also gives
$|x_v'(t)|\geq\frac{|v|}{2}$; moreover, $|x_v(t)|\leq2|v|<\rho$. For fixed $v$, the
state $(x_v,x_v')$ thus remains in a compact set contained in the domain of the
spray on the slit tangent bundle. Theorem~\ref{teo:criterio-prolongacion-edo-banach} rules out
a maximal endpoint before $1$, and the argument reaches all of $[0,1]$. Integrating once
more gives
\begin{equation}\label{eq:estimacion-exponencial-finsler-cero}
 |\widehat\exp_p(v)-v|\leq C_0|v|^2,
 \qquad C_0:=8M_0.
\end{equation}

For the derivative with respect to the initial value, set
$J_v(t):=D_vx_v(t)$. The variational equation is
\begin{equation}\label{eq:variacional-spray-finsler-cero}
 J_v''=-2D_xG(x_v,x_v')J_v-2D_vG(x_v,x_v')J_v',
 \qquad J_v(0)=0,\quad J_v'(0)=I.
\end{equation}
The preceding bounds give
$\|D_xG(x_v,x_v')\|\leq4M_x|v|^2$ and
$\|D_vG(x_v,x_v')\|\leq2M_v|v|$. Applying
Grönwall's lemma~\ref{lema: gronwall} to the first-order
system for $(J_v,J_v')$ and comparing it with $(tI,I)$ gives
\begin{equation}\label{eq:estimacion-diferencial-exponencial-finsler-cero}
 \|D\widehat\exp_p(v)-I\|
 \leq C_1|v|,
 \qquad
 C_1:=(8M_xr+4M_v)
       \exp(1+8M_xr^2+4M_vr).
\end{equation}
In detail, set $E_0:=\exp(1+8M_xr^2+4M_vr)$.
Grönwall's lemma~\ref{lema: gronwall} applied to the sum norm of the system gives
$\|J_v(t)\|+\|J_v'(t)\|\leq E_0$ for $t\in[0,1]$. Integrating
\eqref{eq:variacional-spray-finsler-cero}, we obtain
\[
 \|J_v'(t)-I\|
 \leq\int_0^t\bigl(8M_x|v|^2\|J_v(s)\|
                   +4M_v|v|\|J_v'(s)\|\bigr)\,ds
 \leq (8M_xr+4M_v)E_0|v|t.
\]
Integrating again and using $J_v(0)=0$ yields
\[
 \|J_v(t)-tI\|\leq\frac12(8M_xr+4M_v)E_0|v|t^2.
\]
At $t=1$ this bound implies
\eqref{eq:estimacion-diferencial-exponencial-finsler-cero} with the indicated value of
$C_1$. Estimates
\eqref{eq:estimacion-exponencial-finsler-cero} and
\eqref{eq:estimacion-diferencial-exponencial-finsler-cero} extend the
differential continuously to the origin with value $I$ and prove that
$d(\exp_p)_0=\operatorname{id}_{T_pM}$. In general, $\exp_p$ is not $C^2$ at
$0$, reflecting that $F^2$ and the spray on the slit tangent bundle need not be smooth at
the zero section. The constants $C_0,C_1$ are independent of the direction
$v$ and depend only on $r$ and the bounds $M_0,M_x,M_v$ on the chosen coordinate
neighborhood and unit sphere. The inverse function theorem,
applied to this $C^1$ map, also shows that $\exp_p$ restricts to a
$C^1$ diffeomorphism between neighborhoods of $0$ and $p$; away from $0$ both the
map and its local inverse are smooth.

On the slit tangent bundle, consider the \textbf{pullback bundle}
\[
 \widetilde\pi^*TM
 =\{((x,v),w)\mid v\in T_xM\setminus\{0\},\ w\in T_xM\}.
\]
The fundamental tensor is a metric on this bundle: on the fiber over $(x,v)$
we use $\mathbf{g}_v$. The tautological section is
$\boldsymbol{\mathcal Y}(x,v)=v$.

\begin{theorem}[Chern connection]
\label{teo:conexion-chern-finsler}
Let $(M,F)$ be a finite-dimensional classical Finsler manifold.
\index{Chern connection}
There exists a unique linear connection $\nabla^{\mathrm{Ch}}$ on
$\widetilde\pi^*TM$ whose connection forms
$\displaystyle\boldsymbol{\omega}^i_j=\displaystyle\sum_{k=1}^n\Gamma^i_{jk}(x,v)\,\mathbf{d}x^k$ satisfy
\begin{align}
 \sum_{j=1}^n\mathbf{d}x^j\wedge\boldsymbol{\omega}^i_j&=0,
 \label{eq:chern-torsion-cero}\\
 \mathbf{d}g_{ij}-\sum_{k=1}^n\bigl(g_{kj}\boldsymbol{\omega}^k_i+g_{ik}\boldsymbol{\omega}^k_j\bigr)
 &=2\sum_{k=1}^n C_{ijk}\,\boldsymbol{\delta v}^{\,k}.
 \label{eq:chern-casi-compatible}
\end{align}
Its coefficients are given by
\begin{equation}\label{eq:coeficientes-chern-finsler}
 \Gamma^i_{jk}
 =\frac12\sum_{\ell=1}^n g^{i\ell}
 \left(
 \boldsymbol{\delta}_{x^j}g_{\ell k}
 +\boldsymbol{\delta}_{x^k}g_{j\ell}
 -\boldsymbol{\delta}_{x^\ell}g_{jk}
 \right),
 \qquad
 \boldsymbol{\delta}_{x^j}
 =\boldsymbol{\partial}_{j}-\sum_{m=1}^n N^m_j\boldsymbol{\partial}_{v^m}.
\end{equation}
If $y=\tau(x)$ is a coordinate change, its full transformation law
is
\begin{equation}\label{eq:transformacion-coeficientes-chern-finsler}
 \sum_{\beta,\gamma=1}^n\widetilde\Gamma^\alpha_{\beta\gamma}(y,D\tau(x)v)
 \frac{\partial y^\beta}{\partial x^i}
 \frac{\partial y^\gamma}{\partial x^j}
 =\sum_{k=1}^n\frac{\partial y^\alpha}{\partial x^k}\Gamma^k_{ij}(x,v)
  -\frac{\partial^2y^\alpha}{\partial x^i\partial x^j}.
\end{equation}
Moreover,
\begin{equation}\label{eq:chern-contraccion-spray}
 N^i_j=\sum_{k=1}^n\Gamma^i_{jk}v^k,
 \qquad
 2G^i=\sum_{j,k=1}^n\Gamma^i_{jk}v^jv^k.
\end{equation}
If $F$ comes from a Riemannian metric, $\mathbf{C}=0$ and
$\nabla^{\mathrm{Ch}}$ is the lift of the Levi--Civita connection.
\end{theorem}

\begin{proof}
The first condition is equivalent to
$\Gamma^i_{jk}=\Gamma^i_{kj}$. Since
$\frac{\partial g_{ij}}{\partial v^k}=2C_{ijk}$, the vertical components of
\eqref{eq:chern-casi-compatible} hold automatically. Its
horizontal components are
\[
 \boldsymbol{\delta}_{x^k}g_{ij}
 =\sum_{m=1}^n\bigl(g_{mj}\Gamma^m_{ik}+g_{im}\Gamma^m_{jk}\bigr).
\]
Fix $i,j,k\in\{1,\ldots,n\}$. The three required identities are
\begin{align*}
 \boldsymbol\delta_{x^j}g_{ik}
 &=\sum_{m=1}^n\bigl(g_{mk}\Gamma^m_{ij}+g_{im}\Gamma^m_{kj}\bigr),\\
 \boldsymbol\delta_{x^k}g_{ij}
 &=\sum_{m=1}^n\bigl(g_{mj}\Gamma^m_{ik}+g_{im}\Gamma^m_{jk}\bigr),\\
 \boldsymbol\delta_{x^i}g_{jk}
 &=\sum_{m=1}^n\bigl(g_{mk}\Gamma^m_{ji}+g_{jm}\Gamma^m_{ki}\bigr).
\end{align*}
Adding the first two and subtracting the third, symmetry in the lower
indices cancels the terms involving $g_{mk}$ and $g_{mj}$. What remains is
$\displaystyle2\displaystyle\sum_{m=1}^n g_{im}\Gamma^m_{jk}$. Multiplication by the
inverse matrix of $g$ gives \eqref{eq:coeficientes-chern-finsler}. This proves
local uniqueness. Conversely, the formula is symmetric in $j,k$ and,
on substituting it into the horizontal identity, the two pairs of terms with
opposite signs cancel, leaving $\boldsymbol\delta_{x^k}g_{ij}$.
The vertical components have already been checked; thus the local
connection also exists.

To justify gluing, consider on $\widetilde\pi^*TM$ the fiber-valued
form that sends a tangent vector $Z$ to
$T\widetilde\pi(Z)$. Its components are $\mathbf d x^i$. The first
equation expresses that its exterior covariant derivative is zero, so it is
invariant under coordinate changes. For the second, the vertical
projection determined by the horizontal distribution has components
$\boldsymbol\delta v^{\,i}$. Law
\eqref{eq:transformacion-conexion-no-lineal-finsler} shows that these
components change by the Jacobian of the coordinate change.
The tensor $\mathbf C$, obtained by three derivatives in a fiber where the
change is linear, is covariant of order three. Hence the right-hand side
of \eqref{eq:chern-casi-compatible} and the covariant derivative of $\mathbf g$
are tensorial forms of the same type. Both conditions are preserved under
coordinate changes. The change of basis in the fiber depends only on $x$,
so it introduces no $\mathbf d v$ components into the connection
forms. Local uniqueness gives the gluing.
Finally, differentiating the change of components of a field along
a curve produces $D^2\tau(\dot\gamma,X)$; moving it to the other side
gives the negative sign in
\eqref{eq:transformacion-coeficientes-chern-finsler}. This is the same transformation
rule already proved in \eqref{eq:ley-transformacion-christoffel-hilbert}.

We prove the contraction with the tautological section
$\boldsymbol{\mathcal Y}(x,v)=v$. Write $L=\frac12F^2$. The Euler identities give
\[
 L_{v^i}=\sum_{k=1}^n g_{ik}v^k,
 \qquad
 \frac{\partial g_{ik}}{\partial v^m}=2C_{ikm},
 \qquad \sum_{k=1}^n C_{ikm}v^k=0.
\]
The Euler--Lagrange equation, written using the spray, is
\begin{equation}\label{eq:euler-lagrange-contraida-spray-finsler}
 2\sum_{\ell=1}^n g_{i\ell}G^\ell=\sum_{k=1}^n L_{x^kv^i}v^k-L_{x^i}.
\end{equation}
Differentiating it with respect to $v^j$ gives
\begin{equation}\label{eq:derivada-euler-lagrange-spray-finsler}
 2\sum_{\ell=1}^n g_{i\ell}N^\ell_j
 =L_{x^jv^i}-L_{x^iv^j}
  +\sum_{k=1}^n\frac{\partial g_{ij}}{\partial x^k}v^k
  -4\sum_{m=1}^n C_{ijm}G^m.
\end{equation}
We abbreviate $j\in\{1,\ldots,n\}$ by $\delta_j:=\boldsymbol\delta_{x^j}$.
On the other hand, contracting
\eqref{eq:coeficientes-chern-finsler} in the index $k$ with
$\boldsymbol{\mathcal Y}^k=v^k$ gives
\begin{align*}
 2\sum_{\ell,k=1}^n g_{i\ell}\Gamma^\ell_{jk}v^k
 &=\sum_{k=1}^n\bigl(\delta_jg_{ik}+\delta_kg_{ij}
   -\delta_ig_{jk}\bigr)v^k\\
 &=L_{x^jv^i}-L_{x^iv^j}
   +\sum_{k=1}^n\frac{\partial g_{ij}}{\partial x^k}v^k
   -4\sum_{m=1}^n C_{ijm}G^m.
\end{align*}
The second equality used $\displaystyle\sum_{k=1}^n C_{ikm}v^k=0$ and
$\displaystyle\sum_{k=1}^n N^m_kv^k=2G^m$. Comparison with
\eqref{eq:derivada-euler-lagrange-spray-finsler} and invertibility of
$(g_{ij})$ prove $\displaystyle N^i_j=\displaystyle\sum_{k=1}^n\Gamma^i_{jk}v^k$. The second equality in
\eqref{eq:chern-contraccion-spray} follows by contracting once more and using
$\displaystyle\sum_{j=1}^n N^i_jv^j=2G^i$. In the Riemannian case, $g_{ij}$ does not depend on $v$ and
\eqref{eq:coeficientes-chern-finsler} is the usual Christoffel formula.
\end{proof}

If $\gamma$ is a regular curve and $\mathbf{W}$ is a nonvanishing field along
$\gamma$, the Chern connection defines
\begin{equation}\label{eq:derivada-chern-referencia}
 D_t^{\mathbf{W}}\mathbf{X}
 =\sum_{i=1}^n\left(\dot X^i+\sum_{j,k=1}^n\Gamma^i_{jk}(\gamma,\mathbf{W})\dot\gamma^jX^k\right)
   \boldsymbol{\partial}_{i}.
\end{equation}
Almost compatibility is expressed by
\begin{equation}\label{eq:casi-compatibilidad-curva-chern}
 \frac d{dt}\mathbf{g}_{\mathbf{W}}(\mathbf{X},\mathbf{Y})
 =\mathbf{g}_{\mathbf{W}}(D_t^{\mathbf{W}}\mathbf{X},\mathbf{Y})
  +\mathbf{g}_{\mathbf{W}}(\mathbf{X},D_t^{\mathbf{W}}\mathbf{Y})
  +2\mathbf{C}_{\mathbf{W}}(\mathbf{X},\mathbf{Y},D_t^{\mathbf{W}}\mathbf{W}).
\end{equation}
In particular, if $\mathbf{T}_\gamma:=\dot\gamma$ is the tangent field, a
regular curve is a geodesic if and only if
$D_t^{\mathbf{T}_\gamma}\mathbf{T}_\gamma=0$; by
\eqref{eq:chern-contraccion-spray}, this is precisely
\eqref{eq:geodesica-spray-finsler}.

\subsection{First variation and directed distance}

For a regular curve, set
\[
 E_F(\gamma):=\frac12\int_a^bF(\dot\gamma)^2dt,
 \qquad
 L_F(\gamma):=\int_a^bF(\dot\gamma)dt.
\]

\begin{proposition}[Finsler first variation]\label{prop:primera-variacion-finsler}
Let $(M,F)$ be a finite-dimensional classical Finsler manifold.
Let $\alpha\colon (-\varepsilon,\varepsilon)\times[a,b]\longrightarrow M$ be a smooth variation
through regular curves, let $\mathbf{T}=\partial_t\alpha$, $\mathbf{V}=\partial_s\alpha$, and let
$\gamma=\alpha_0$. On the right-hand sides, $\mathbf{T}$ and $\mathbf{V}$ are evaluated at $s=0$.
Then
\[
 \frac d{ds}\biggr|_{0}\frac12\int_a^bF(\dot\alpha_s)^2dt
 =\left[\sum_{i=1}^n\frac{\partial L}{\partial v^i}V^i\right]_a^b
 +\int_a^b\sum_{i=1}^n\left(
 \frac{\partial L}{\partial x^i}
 -\frac d{dt}\frac{\partial L}{\partial v^i}
 \right)V^i dt.
\]
Equivalently,
\begin{equation}\label{eq:primera-variacion-energia-chern}
 \frac d{ds}\biggr|_{0}E_F(\alpha_s)
 =\bigl[\mathbf{g}_{\mathbf{T}}(\mathbf{V},\mathbf{T})\bigr]_a^b
 -\int_a^b \mathbf{g}_{\mathbf{T}}(\mathbf{V},D_t^{\mathbf{T}}\mathbf{T})\,dt.
\end{equation}
For length we have
\begin{equation}\label{eq:primera-variacion-longitud-finsler}
 \frac d{ds}\biggr|_{0}L_F(\alpha_s)
 =\left[\frac{\mathbf{g}_{\mathbf{T}}(\mathbf{V},\mathbf{T})}{F(\mathbf{T})}\right]_a^b
 -\int_a^b \mathbf{g}_{\mathbf{T}}\left(\mathbf{V},D_t^{\mathbf{T}}\frac{\mathbf{T}}{F(\mathbf{T})}\right)dt.
\end{equation}
Thus, the energy extremals with fixed endpoints are geodesics with an
affine parameter. The regular length extremals are curves whose arc length
parametrization is a geodesic.
\end{proposition}

\begin{proof}
Fix $0<\delta<\varepsilon$ and restrict the variation to
$[-\delta,\delta]\times[a,b]$. Since it is regular and this rectangle is compact,
the map $(s,t)\mapsto(\alpha(s,t),\mathbf{T}(s,t))$ has compact image in
$\widetilde{TM}$; the derivatives of $L=\frac12F^2$ are continuous and bounded on
that image. Cover the central curve by finitely many charts and subdivide
$[a,b]$ so that, after shrinking $\varepsilon$, each strip of the
variation lies in one of them. Differentiation under the integral sign is
therefore justified and gives on each subinterval
\[
 \frac d{ds}\biggr|_0E_F(\alpha_s)
 =\int_a^b\sum_{i=1}^n\left(
 \frac{\partial L}{\partial x^i}V^i
 +\frac{\partial L}{\partial v^i}\dot V^i
 \right)dt.
\]
Integration by parts is written explicitly as
\[
 \int_a^b\sum_{i=1}^n\frac{\partial L}{\partial v^i}\dot V^i\,dt
 =\left[\sum_{i=1}^n\frac{\partial L}{\partial v^i}V^i\right]_a^b
  -\int_a^b\sum_{i=1}^n\frac d{dt}\left(\frac{\partial L}{\partial v^i}\right)V^i\,dt,
\]
and proves the first formula on summing over the subintervals: the interior boundary
terms cancel because
$\displaystyle\sum_{i=1}^n L_{v^i}V^i=d_vL(\mathbf{V})$ is the intrinsic evaluation of a covector. Now
$\displaystyle L_{v^i}=\displaystyle\sum_{j=1}^n g_{ij}v^j$ and
\eqref{eq:euler-lagrange-contraida-spray-finsler} imply
\[
 \frac d{dt}L_{v^i}-L_{x^i}
 =\sum_{j=1}^n g_{ij}\bigl(\ddot x^j+2G^j(x,\dot x)\bigr)
 =\sum_{j=1}^n g_{ij}(D_t^{\mathbf{T}}\mathbf{T})^j.
\]
The last equality uses
\eqref{eq:chern-contraccion-spray}. Since the boundary term is
$\displaystyle\sum_{i=1}^n L_{v^i}V^i=\mathbf{g}_{\mathbf{T}}(\mathbf{T},\mathbf{V})$, this proves
\eqref{eq:primera-variacion-energia-chern}.

For length, the Euler identity in
Proposition~\ref{prop:euler-legendre-finsler} gives
$d(F_x)_{\mathbf{T}}(\mathbf{W})=\frac{\mathbf{g}_{\mathbf{T}}(\mathbf{T},\mathbf{W})}{F(\mathbf{T})}$.
Symmetry of the lower indices of
$\Gamma^i_{jk}$, together with
$\partial_s\mathbf{T}=\partial_t\mathbf{V}$, produces the zero-torsion
identity along the variation
\[
 D_s^{\mathbf{T}}\mathbf{T}=D_t^{\mathbf{T}}\mathbf{V}.
\]
Applying \eqref{eq:casi-compatibilidad-curva-chern} with reference $\mathbf{T}$ eliminates the
Cartan term because one of its arguments is $\mathbf{T}$; therefore,
\[
 \partial_sF(\mathbf{T})
 =\mathbf{g}_{\mathbf{T}}\left(D_t^{\mathbf{T}}\mathbf{V},\frac{\mathbf{T}}{F(\mathbf{T})}\right).
\]
The same almost compatibility formula, now in the direction $t$, gives
\[
 \frac d{dt}\mathbf{g}_{\mathbf{T}}\left(\mathbf{V},\frac{\mathbf{T}}{F(\mathbf{T})}\right)
 =\mathbf{g}_{\mathbf{T}}\left(D_t^{\mathbf{T}}\mathbf{V},\frac{\mathbf{T}}{F(\mathbf{T})}\right)
  +\mathbf{g}_{\mathbf{T}}\left(\mathbf{V},D_t^{\mathbf{T}}\frac{\mathbf{T}}{F(\mathbf{T})}\right),
\]
since $\mathbf{C}_{\mathbf{T}}\left(\mathbf{V},\frac{\mathbf{T}}{F(\mathbf{T})},
D_t^{\mathbf{T}}\mathbf{T}\right)=0$. Integrating this identity proves
\eqref{eq:primera-variacion-longitud-finsler}, with the indicated boundary
terms.

It remains to justify that the fundamental lemma applies to all the required
local fields. Let $J\subset(a,b)$ be a compact subinterval contained in a
chart $\varphi\colon U\longrightarrow\mathbb R^n$ in the sense that
$\gamma(J)\subseteq U$, and let $\mathbf{W}$ be a smooth field
along $\gamma$ supported in the interior of $J$. For small $|s|$,
\[
 \alpha(s,t):=\varphi^{-1}\!\left(
  \varphi(\gamma(t))+s\,d\varphi_{\gamma(t)}\mathbf{W}(t)\right)
\]
on $J$, extended by $\gamma$ outside $J$, is a smooth variation.
Compactness of the support allows a uniform choice of $s$, the gluing is
smooth because $\mathbf{W}$ vanishes near the endpoints of $J$, and regularity of
$\gamma$ ensures that nearby curves remain regular after shrinking
$|s|$.

For fixed endpoints the boundary terms disappear. Fixing
$z\in\mathbb R^n$ and $\chi\in C_c^\infty(\mathring J)$, choose
$\mathbf{W}(t)=\chi(t)(d\varphi_{\gamma(t)})^{-1}z$. The energy formula and the
fundamental lemma imply
\[
 \mathbf{g}_{\mathbf{T}}\bigl((d\varphi_{\gamma(t)})^{-1}z,
 D_t^{\mathbf{T}}\mathbf{T}\bigr)=0
 \qquad\text{on }\mathring J.
\]
The arbitrariness of $z$ and positivity of $\mathbf{g}_{\mathbf{T}}$ give
$D_t^{\mathbf{T}}\mathbf{T}=0$, which by
\eqref{eq:chern-contraccion-spray} is
\eqref{eq:geodesica-spray-finsler}. Applied to the length formula, the same
argument gives $D_t^{\mathbf{T}}\left(\frac{\mathbf{T}}{F(\mathbf{T})}\right)=0$;
after arc length parametrization, this is
the affine geodesic equation. The converses follow directly from the
formulas.

Metric differentiation and integration by parts agree with the Levi--Civita
calculation of the finite-dimensional chapter, whose connection and derivative
along curves are given in
Theorem~\ref{teo:variedades-riemannianas-teorema-fundamental-de-la-geometria-riemanniana}
and Theorem~\ref{teo: derivada covariante a lo largo de una curva}. Here we
also use the identity $\mathbf C_{\mathbf T}(\mathbf T,\cdot,\cdot)=0$:
it eliminates the additional terms in Chern almost compatibility and
allows the first variation to be expressed through covariant acceleration.
\end{proof}

The \textbf{forward directed distance} is
\begin{equation}\label{eq:distancia-dirigida-finsler}
 d_F^+(p,q):=\inf\{L_F(\gamma)\mid
 \gamma\text{ is a piecewise $C^1$ curve from $p$ to $q$}\}.
\end{equation}
It satisfies $d_F^+(p,q)>0$ for $p\neq q$ and the triangle inequality, but is not
necessarily symmetric. The reverse metric
\[
 \overleftarrow F(v):=F(-v)
\]
satisfies
$d_{\overleftarrow F}^+(p,q)=d_F^+(q,p)$. We denote this last distance
by $d_F^-(p,q)$. The functions
\[
 d_F^{\max}(p,q):=\max\{d_F^+(p,q),d_F^-(p,q)\},
 \qquad
 d_F^{\mathrm{sum}}(p,q):=d_F^+(p,q)+d_F^-(p,q)
\]
are finite symmetric distances on each connected component. For
points in different components, we adopt the value $+\infty$;
the metric assertions apply componentwise.

\begin{proposition}[Topology of the directed distance]
\label{prop:topologia-distancia-dirigida-finsler}
Let $(M,F)$ be a finite-dimensional classical Finsler manifold.
Forward directed balls, backward directed balls, and
balls of $d_F^{\displaystyle\max}$ and $d_F^{\mathrm{sum}}$ induce the topology of $M$.
\end{proposition}

\begin{proof}
Fix a Euclidean norm $|\,\cdot\,|$ in a chart about $x_0$.
By compactness of the unit sphere and continuity of $F$, after shrinking the
chart there exist $0<m\leq M$ such that
\[
 m|v|\leq F_x(v)\leq M|v|,
 \qquad
 m|v|\leq F_x(-v)\leq M|v|
\]
for every $x$ in the chart and every $v$. The first exit argument in the
proof of Theorem~\ref{teo:distancia-banach-finsler-topologia}, applied
separately to $F$ and $\overleftarrow F$, gives a coordinate ball
$V$ such that, for $x,y\in V$,
\[
 m|\varphi(y)-\varphi(x)|
 \leq d_F^+(x,y)\leq M|\varphi(y)-\varphi(x)|
\]
and the same inequalities with $d_F^-$ in place of $d_F^+$. Finally,
\[
 d_F^{\max}(x,y)
 \leq d_F^{\mathrm{sum}}(x,y)
 \leq2d_F^{\max}(x,y).
\]
These estimates prove all four topological assertions.
\end{proof}

\begin{definition}[Directed and geodesic completeness]
\label{def:completitudes-dirigidas-finsler}
\index{completeness!forward Finsler}
Let $(M,F)$ be a classical Finsler manifold. A sequence $(x_n)_{n\in\mathbb N}$ is \textbf{forward Cauchy} if, for every
$\varepsilon>0$, there exists $N\in\mathbb N$ such that
\[
 d_F^+(x_i,x_j)<\varepsilon
 \qquad\text{if }i,j\in\mathbb N,\ j\geq i\geq N.
\]
The manifold is \textbf{forward metrically complete} if every forward
Cauchy sequence converges. Backward metric completeness is defined by
replacing $F$ by $\overleftarrow F$. It is
\textbf{forward geodesically complete} if every maximal geodesic with
nonzero initial velocity is defined for all positive times; the
backward notion is defined analogously.
\end{definition}

Convergence in the preceding definition can be expressed in the topology of
$M$, by $d_F^+(x_n,x)\to0$, or by $d_F^-(x,x_n)\to0$. Indeed, the last two
quantities are equal by definition, and in a smaller chart the
preceding proof gives constants $m,M>0$ such that
\[
 m|\varphi(x_n)-\varphi(x)|
 \leq d_F^+(x_n,x)
 \leq M|\varphi(x_n)-\varphi(x)|
\]
when $x_n$ and $x$ belong to that chart. Forward completeness does not
in general imply backward completeness. If $F$ is reversible, both
distances agree and ordinary metric completeness is recovered.

\begin{example}[Randers metrics]
If $a$ is Riemannian and $\beta$ is a $1$-form with
$|\beta|_a<1$, then
\[
 F(v)=\sqrt{a(v,v)}+\beta(v)
\]
is classical Finsler. It is reversible if and only if $\beta=0$. The condition
$|\beta|_a<1$ is precisely pointwise positivity preventing $F$
from vanishing in a nonzero direction. Over a noncompact base, a uniform
bound requires the additional condition $\|\beta\|_{L^\infty(M,a)}<1$.
To check vertical convexity, fix a fiber and write
$\alpha(v)=\sqrt{a(v,v)}$. If $v\ne0$, then
\[
 \frac12D^2(F^2)_v(w,w)
 =\bigl(d\alpha_v(w)+\beta(w)\bigr)^2
  +\frac{F(v)}{\alpha(v)}
      \left(a(w,w)-\frac{a(v,w)^2}{a(v,v)}\right).
\]
The second summand is nonnegative by Cauchy--Schwarz and vanishes only
when $w$ is proportional to $v$. In that case, if $w=tv\ne0$, the first
summand is $t^2F(v)^2>0$. If they are not proportional, the second summand is
positive. Thus, the Hessian is positive definite for every $v\ne0$.
Moreover, if $b(x):=|\beta_x|_a<1$, we have
$(1-b(x))\alpha(v)\leq F(v)\leq(1+b(x))\alpha(v)$; this distinguishes
pointwise positivity from uniform comparison over the whole base.
\end{example}

\begin{example}[Funk metric]
On the Euclidean ball $B^n$, the formula
\[
 F_x(v)=\frac{\sqrt{(1-|x|^2)|v|^2+\langle x,v\rangle^2}
                 +\langle x,v\rangle}{1-|x|^2}
\]
defines a nonreversible metric. It is a Randers metric: its Riemannian
part and its $1$-form are
\[
 a_x(v,v)=\frac{|v|^2}{1-|x|^2}
          +\frac{\langle x,v\rangle^2}{(1-|x|^2)^2},
 \qquad \beta_x(v)=\frac{\langle x,v\rangle}{1-|x|^2}.
\]
The inverse matrix of $a_x$ is $(1-|x|^2)(I-xx^{\mathsf T})$, so
$|\beta_x|_{a_x}^2=|x|^2<1$. The preceding example ensures strong vertical
convexity. Along a radius directed toward the boundary,
$F_{te_1}(e_1)=\frac{1}{1-t}$, whose integral diverges; in the opposite direction one obtains
$F_{te_1}(-e_1)=\frac{1}{1+t}$. More generally, if the ray starting at $x$ and passing
through $y$ meets the boundary at $b$, we have
\[
 d_F^+(x,y)=\log\frac{|x-b|}{|y-b|}.
\]
To verify the formula, fix $x\ne y$ and that point $b\in\partial B^n$,
and set $h_b(z)=1-\langle b,z\rangle>0$. For $w\ne0$, the formula for
$F_z(w)$ is the positive solution of
$|z+w/F_z(w)|=1$. Taking the inner product with $b$ gives
\[
 F_z(w)\geq\frac{\langle b,w\rangle}{h_b(z)}
 =-d(\log h_b)_z(w).
\]
The inequality also holds for $w=0$. Every curve from $x$ to $y$ therefore has
length at least $\log(h_b(x)/h_b(y))$. On the oriented segment
from $x$ to $y$ equality holds, because the tangent ray
meets the sphere precisely at $b$. Since $x,y,b$ are collinear,
$h_b(x)/h_b(y)=|x-b|/|y-b|$, proving the distance formula.
In particular,
\[
 d_F^+(0,y)=-\log(1-|y|),
 \qquad \overline B^+(0,R)=\{y\in B^n\mid |y|\leq1-e^{-R}\},
 \qquad R\geq0.
\]
These balls are compact. A forward Cauchy sequence lies
in one of them by the triangle inequality and has a
subsequence converging to $y\in B^n$. To check convergence of
the whole sequence, fix $\varepsilon>0$ and an index $N$ beyond which
$d_F^+(x_i,x_j)<\varepsilon/2$ if $j\geq i\geq N$. For each $i\geq N$,
take subsequence indices greater than $i$ and pass to the limit;
continuity of the distance gives $d_F^+(x_i,y)\leq\varepsilon/2$.
The topology of backward balls implies $x_i\to y$. This gives
forward completeness. For the reverse metric,
however, the radius $t\mapsto te_1$ has length
\[
 \int_0^1F_{te_1}(-e_1)dt=\log2<\infty.
\]
Thus, a sequence along this radius converging to the boundary is Cauchy for
the reverse distance and does not converge in $B^n$; the metric is not backward
complete.
\end{example}

\begin{theorem}[Finite-dimensional Finsler Hopf--Rinow theorem]
\label{teo:hopf-rinow-finsler-clasico}
Let $(M,F)$ be a connected finite-dimensional classical Finsler manifold. The following are
equivalent:
\begin{enumerate}[label=\textup{(\roman*)}]
 \item $(M,F)$ is forward metrically complete;
 \item $(M,F)$ is forward geodesically complete;
 \item $\exp_p$ is defined on all of $T_pM$ for each $p\in M$;
 \item every forward closed ball
       $\overline B^+(p,r)=\{q\in M\mid d_F^+(p,q)\leq r\}$ is compact;
 \item every forward bounded closed subset, that is, one contained
       in some ball $B^+(p,r)$, is compact.
\end{enumerate}
Under these conditions, for any $p,q\in M$ there exists a directed geodesic
from $p$ to $q$ whose length is $d_F^+(p,q)$. The analogous assertions
hold backward, applied to $\overleftarrow F$. If $F$ is
reversible, the two versions agree and give the same completeness
equivalences as the Riemannian Hopf--Rinow theorem;
a reversible metric need not come from an inner product.
\end{theorem}

A proof can be found in
\cite[\S 6.6, Thm.~6.6.1]{BaoChernShen2000}; it combines local existence of
geodesics with compactness arguments for minimizing sequences and properties of
directed balls. Finite dimensionality is essential both for
compactness of indicatrices and for the local compactness used in
the argument, and cannot be omitted in view of
Remark~\ref{obs:hopf-rinow-infinito}. The result relates, in each orientation separately, metric
completeness, geodesic completeness, and existence of minimizing geodesics.
It does not identify forward completeness with backward completeness.

\section{Global results and distinctions between theories}
\label{sec:comparacion-riemann-finsler-banach}

A smooth spray on all of $TM$ determines the symmetric connection used in
the following result. In coordinates we write, as in
Definition~\ref{def:spray-banach},
\[
 \mathbf{S}(x,v)=\bigl((x,v),(v,Q_x(v))\bigr).
\]
Smoothness at the zero section and homogeneity for all scalars
imply that $Q_x$ is quadratic. Its polarization defines
\begin{equation}\label{eq:conexion-asociada-spray-neeb}
 \Gamma_x(u,v)
 :=-\frac12\bigl(Q_x(u+v)-Q_x(u)-Q_x(v)\bigr),
 \qquad Q_x(v)=-\Gamma_x(v,v).
\end{equation}
The transformation law for the spray gives
\eqref{eq:ley-transformacion-christoffel-hilbert} for these coefficients, without
using an inner product; hence they define a symmetric linear connection on
$TM$. Its parallel transport is that of
Theorem~\ref{teo:transporte-paralelo-haz-banach}. We say that the spray is
\textbf{compatible} with a Banach--Finsler structure $F$ if this transport,
along each geodesic segment, preserves the norms $F_x$.

\begin{theorem}[Neeb's Cartan--Hadamard theorem]
\label{teo:cartan-hadamard-neeb}
Let $M$ be a connected smooth Banach manifold equipped with a
Banach--Finsler structure $F$ and a smooth spray on all of $TM$, compatible
with $F$ in the preceding sense. Suppose it has nonpositive
curvature in the following sense: for every $x\in M$ and every
$v\in\mathcal D_x$, where $\mathcal D_x$ is the domain of $\exp_x$, and every
$w\in T_xM$,
\[
 (d\exp_x)_v\colon T_xM\longrightarrow T_{\exp_x(v)}M
 \quad\text{is an isomorphism and}\quad
 F_{\exp_x(v)}((d\exp_x)_vw)\geq F_x(w).
\]
Under these hypotheses, geodesic completeness is equivalent to completeness of
$(M,d_F)$. If either holds, then for every $x\in M$ the map
$\exp_x\colon T_xM\longrightarrow M$ is a covering map. If $M$ is also
simply connected,
$\exp_x$ is a global diffeomorphism. In this last case radial geodesics
minimize distance between their endpoints.
\end{theorem}

The notions and results that follow come from Neeb~\cite[pp.~118--124]{Neeb2002}. Compatible tangent norms and compatible sprays are introduced in his Definitions~1.3--1.4, pp.~118--120. The implication from metric
completeness to geodesic completeness is his Lemma~1.8, and the path lifting
argument is his Theorem~1.9, pp.~121--123. The covering property of the exponential
map is Theorem~1.10, and equivalence of the completeness notions is
Corollary~1.11, pp.~123--124; the diffeomorphism and
minimization in the simply connected case appear in Corollary~1.12, p.~124. Neeb's argument constructs
path lifts with length control, extends geodesics, and
verifies the covering map criterion. The estimate making that
argument possible is
\[
 F_x\bigl(((d\exp_x)_v)^{-1}z\bigr)\leq F_{\exp_x(v)}(z),
 \qquad z\in T_{\exp_x(v)}M.
\]
Thus, the length of a lift does not exceed that of its projection;
in the complete domain $T_xM$, a lift of finite length has
a limit on reaching a finite endpoint, and a local inverse allows it to be
continued. This is the path continuation mechanism used in
Chapter~\ref{cap:calculo-diferencial-banach}. The condition on $d\exp_x$ contains two
independent assertions, invertibility and expansiveness, neither of which
follows from continuity of $F$ alone. In turn, a continuous Banach--Finsler
structure does not by itself determine a compatible spray.

\begin{proof}
Suppose first that $(M,d_F)$ is complete. If a maximal geodesic
$\gamma$ had a finite right endpoint $b$, compatibility of the spray
would make its speed $F_{\gamma(t)}(\dot\gamma(t))=a$ constant. Indeed,
its velocity is parallel for the connection associated with the spray, and transport
along the geodesic itself preserves $F$. Hence,
$d_F(\gamma(s),\gamma(t))\leq a|t-s|$, and completeness provides a
limit $q$ as $t\to b^{-}$.

In a chart about $q$, local comparison of norms gives
$c>0$ such that $c\|w\|_X\leq F_z((d\varphi_z)^{-1}w)$.
Smoothness of the connection coefficients allows the chart to be shrunk
so that $\|\Gamma_z\|\leq K$, with $K>0$. If
$u(t)=\varphi(\gamma(t))$, on the tail of the curve we have
\[
 \|u'(t)\|_X\leq a/c,
 \qquad \|u''(t)\|_X\leq K(a/c)^2.
\]
Hence $u'$ has a limit in $X$ as $t\to b^{-}$. The state
$(u,u')$ converges to a point of the spray domain; the ODE continuation
criterion in the appendix extends the geodesic. This contradicts
maximality. The same proof at the left endpoint gives geodesic
completeness. This is the argument of
Proposition~\ref{prop:completo-metrico-implica-geodesico-hilbert}, with
constant speed and Finsler comparison in place of Riemannian
compatibility.

Now suppose geodesic completeness, fix $x\in M$, and write
$E:=T_xM$, with the complete norm $F_x$, and $f:=\exp_x\colon E\to M$.
The differential of $f$ is an isomorphism at every point, so the
inverse function theorem makes it a local diffeomorphism. If
$c\colon[0,1]\to M$ is a piecewise $C^1$ curve and $v_0\in E$ satisfies
$f(v_0)=c(0)$, a local inverse starts a lift $v$ of $c$ with
$v(0)=v_0$. Two lifts with the same value agree near that
time by local injectivity of $f$; the open-and-closed argument
extends agreement to their common interval.

On every segment where the lift is defined, the chain rule
and expansiveness give
\[
 F_x(v'(t))
 =F_x\bigl((df_{v(t)})^{-1}c'(t)\bigr)
 \leq F_{c(t)}(c'(t)).
\]
If its initial domain had an endpoint $b\leq1$, then, for
$s<t<b$,
\[
 F_x(v(t)-v(s))
 \leq\int_s^t F_{c(\tau)}(c'(\tau))\,d\tau.
\]
The right-hand side tends to zero as $s,t\to b^{-}$. Completeness
of $E$ provides a limit $v_b$, and continuity of $f$ gives
$f(v_b)=c(b)$. A local inverse about $v_b$ continues the
lift, or defines it at $t=1$ if $b=1$. Thus, every piecewise $C^1$
curve has a unique lift on all of $[0,1]$ with the prescribed initial
value. Since connectedness of $M$ allows $x=f(0)$ to be joined to any point
by a curve of this class, $f$ is surjective.

We verify the covering map property. Let $q\in M$ and choose a
chart $\psi\colon V\to B_X(0,r)$ with $\psi(q)=0$ whose image is an
open ball. For $y\in V$, let
$c_y(t):=\psi^{-1}(t\psi(y))$. For each $e\in f^{-1}(q)$, define
$s_e(y)$ to be the endpoint of the lift of $c_y$ starting at $e$.
Then $f\circ s_e=\operatorname{id}_V$.
These maps are smooth. To check this near a point $y_0\in V$,
cover the compact lift of $c_{y_0}$ by finitely many
domains of local inverses of $f$. Subdivide $[0,1]$ so that each
segment lifts into one of them. On the first segment use the inverse
that takes $q$ to $e$. If the first $j$ segments have been constructed for
$y$ near $y_0$, continuity of the endpoint and openness of the next
domain allow that neighborhood of $y_0$ to be shrunk and the next
local inverse to be used. Joint continuity of $(t,y)\mapsto c_y(t)$ and
compactness of each interval of the subdivision allow this choice to be made
for the entire segment at once. Uniqueness makes the two expressions agree at their common
time. Repeating this step finitely many times expresses the endpoint
through smooth local inverses. This proves smoothness of $s_e$.

Each $s_e(V)$ is open. Indeed, about $s_e(y)$ take a
domain $O$ on which $f$ is a diffeomorphism; by continuity of $s_e$,
there is an open set $W\subseteq V$ about $y$ with $s_e(W)\subseteq O$.
There $s_e=(f\restriction_O)^{-1}$, so $s_e(W)$ is open.
If $s_e(y)=s_{e'}(y)$, lifting the reverse curve of $c_y$ from that
point, uniqueness gives $e=e'$. Hence the images of two distinct
sections are disjoint. Finally, given $z\in f^{-1}(V)$, lift
from $z$ the reverse curve of $c_{f(z)}$; its endpoint is some
$e\in f^{-1}(q)$, and uniqueness gives $z=s_e(f(z))$. We have proved
\[
 f^{-1}(V)=\coprod_{e\in f^{-1}(q)}s_e(V),
 \qquad f\restriction_{s_e(V)}\colon s_e(V)\to V
 \text{ is a diffeomorphism}.
\]
This is the covering map condition.

The same lifting estimate proves metric completeness.
Let $(q_j)_{j\in\mathbb N}$ be a Cauchy sequence in $M$. Choose a
subsequence $(q_{j_k})_{k\in\mathbb N}$ such that
$d_F(q_{j_k},q_{j_{k+1}})<2^{-k-1}$. By the definition of distance,
there exists a curve from $q_{j_k}$ to $q_{j_{k+1}}$ of length less than
$2^{-k}$. Take $v_1\in f^{-1}(q_{j_1})$ and lift these curves
successively, using the endpoint of each as the initial value of the
next. If $v_k$ is the lift of $q_{j_k}$, then
\[
 F_x(v_\ell-v_k)
 \leq\sum_{m=k}^{\ell-1}2^{-m}\leq2^{1-k},
 \qquad k,\ell\in\mathbb N,\quad\ell>k.
\]
Thus, $v_k$ converges in $E$ to a point $v$. Continuity of $f$ and
agreement of the topology of $M$ with that of $d_F$ imply
$q_{j_k}\to f(v)$ in distance. A Cauchy sequence with a
convergent subsequence converges to the same limit, so $(M,d_F)$
is complete.

If $M$ is simply connected, the preceding covering has a single
sheet, because $E$ is connected. This can be seen through lifts: a path
in $E$ between two points in the same fiber projects to a
contractible loop. Lifting a contraction with fixed endpoints makes the lifted
endpoint vary continuously in the discrete fiber, so it remains
constant. The lift of the constant loop ends where it begins;
the two points of the fiber must coincide. The lift of that
contraction is constructed in the evenly covered neighborhoods in the
same way as for paths, subdividing the compact parameter square
into finitely many rectangles. Thus, $f$ is bijective and its
local inverses form a global smooth inverse.

It remains to prove minimization. In the simply connected case, let
$\Gamma=\exp_x^{-1}\circ\gamma$ be the lift
of a curve from $x$ to $y$; then $\Gamma(0)=0$. Expansiveness gives
\[
 L_F(\gamma)\geq\int_0^1F_x(\dot\Gamma(t))dt
 \geq F_x(\Gamma(1)).
\]
The radial curve $t\mapsto\exp_x(t\Gamma(1))$ is a geodesic, and its velocity
is parallel; compatibility of the spray with $F$ shows that its length is
$F_x(\Gamma(1))$. Therefore, it is minimizing.
\end{proof}

\begin{remark}[Four levels that must be distinguished]
\label{obs:cuatro-niveles-geometria-infinita}
\begin{enumerate}
 \item A strong Riemannian metric exists only on Hilbertizable models and
 provides $\sharp_{\mathbf{g}}$, Levi--Civita, and a smooth spray.
 \item A strictly weak Riemannian metric may be geometrically
 natural, such as $L^2$ on mapping spaces of positive regularity. Its musical
 map fails to be surjective on at least one fiber; being weak, in this
 chapter's convention, does not rule out being strong.
 \item A Palais--Struwe Banach--Finsler structure is a locally uniform family
 of norms; it requires neither a vertical Hessian nor a Chern
 connection, and is the appropriate structure for completeness and deformations.
 \item A classical Finsler metric requires smoothness away from zero and strong
 vertical convexity. The results developed here also use
 compactness of spheres and finite dimensionality of the vertical tensors;
 their Banach analogues require additional smoothness and
 uniform convexity hypotheses and do not follow from an arbitrary norm.
\end{enumerate}
No tensor product, musical duality, or identification of a space with
its dual is regarded as canonical without first fixing the structure that produces it.
\end{remark}

\chapter{Uniform geometry and geometric function spaces}
\label{cap:geometria-uniforme-funcional}

The preceding chapters established differential calculus and the basic geometric structures on Banach manifolds. To apply them to spaces of functions and sections, an additional difficulty remains: charts, metrics, and connections must be compatible with the chosen functional regularity, and their estimates must remain controlled as the basepoint varies.

In a mapping space, a chart is constructed pointwise using a local addition or the exponential map of the target space. In spaces of metrics and connections, the structure comes from open subsets or affine spaces of tensor sections. Gauge and diffeomorphism groups also introduce nonlinear operations with different regularity: the gauge action requires one derivative of the group element, whereas diffeomorphism composition loses derivatives in one variable. These differences determine which objects form Hilbert manifolds at fixed regularity and which must be studied through a Sobolev scale or an ILH structure.

The aim of this chapter is to bring these constructions together in a uniform framework. We will identify the charts, tangent spaces, and natural actions, while specifying in each case the regularity assumptions that make the operations well defined and smooth. This precision will be essential in the next chapter: many geometric energies use an \(L^2\) metric that is weak relative to the Sobolev topology, and their differential gradients lower regularity. Infinite-dimensional geometry must therefore retain both the geometric structure and the analytic scale in which the equations live.

\begin{semblanzaHistorica}{Eichhorn, Amann, and control at infinity}
Noncompact geometry requires distinguishing local closeness, uniform control, and global distance. Eichhorn developed uniformities for comparing geometric objects and studying their finite-distance components; Amann organized uniformly regular atlases adapted to the analysis of function spaces. Although both languages share uniform estimates, they address different questions: one describes configuration spaces\glsadd{espacio-configuraciones}, while the other allows operators and norms to be localized without the constants deteriorating as one escapes to infinity.
\end{semblanzaHistorica}

\section{Uniform control of the base manifold: uniform regularity and bounded geometry}
\label{sec:control-uniforme-base-amann}

The spaces $BC^j(U,F)$ and $BUC^j(U,F)$, including the conventions for a
half-box and for $j=0,\infty$, were defined in
Definition~\ref{def:espacios-bc-buc}. Normalized and shrinkable atlases, their
multiplicity, and uniformly regular atlases and structures were
defined in Definitions~\ref{def:atlas-uniformemente-regular-amann}
and~\ref{def:estructura-uniformemente-regular-amann}. For a chart
$\kappa\colon U_\kappa\longrightarrow Q_\kappa$, the coordinate representation
of the metric is, by definition,
$g_\kappa=(\kappa^{-1})^*g$. The uniform ellipticity constants and
$BUC^j$ bounds for $g_\kappa$ and $g_\kappa^{-1}$ appear explicitly in
\eqref{eq:amann-comparacion-metrica-explicita},
\eqref{eq:amann-buc-metrica}, and
\eqref{eq:amann-buc-metrica-inversa}; none depends on the chart.

Theorem~\ref{teo:geometria-acotada-implica-regularidad-uniforme-amann}
and
Corollary~\ref{cor:geometria-acotada-frontera-regularidad-uniforme-amann}
prove the direction used in this chapter, without and with boundary. In the
first case, the constants depend only on the dimension, a common normal
radius, and bounds on curvature derivatives of the required
order. In the second, one adds the collar width, the normal radius of the
boundary, and bounds on the second fundamental form and its derivatives. The
full comparisons and precise references for their converses are
stated in
Theorems~\ref{teo:amann-ur-equivale-geometria-acotada}
and~\ref{teo:amann-ur-equivale-geometria-acotada-frontera}, respectively.

Amann's localization system was constructed in
Proposition~\ref{prop:sistema-localizacion-amann}. Without boundary, for the
analysis operators $\mathcal S_{\mathfrak K}$ and synthesis operators
$\mathcal R_{\mathfrak K}$ defined in
Corollary~\ref{cor:puente-funcional-amann-sobolev}, we have
$\mathcal R_{\mathfrak K}\mathcal S_{\mathfrak K}=\operatorname{id}$ and,
for each $s\in\mathbb R$ and $1<p<\infty$, there exist constants
$c_{s,p},C_{s,p},D_{s,p}>0$ such that
\[
c_{s,p}\|u\|_{H^{s,p}(M)}
\leq
\|\mathcal S_{\mathfrak K}u\|_{\ell^p(I;H^{s,p}(\mathbb R^m))}
\leq C_{s,p}\|u\|_{H^{s,p}(M)},
\]
as well as
$\|\mathcal R_{\mathfrak K}v\|_{H^{s,p}(M)}
\leq D_{s,p}\|v\|_{\ell^p(I;H^{s,p}(\mathbb R^m))}$.
The three constants depend only on $m,s,p$, the multiplicity, and
the uniform bounds on coordinate changes and cutoff functions through the
required orders; they do not depend on the chart or on $u$ or $v$. The same
corollary records Triebel--Lizorkin comparisons and, when the
indices allow, Besov comparisons with explicit constants.

With boundary,
Corollary~\ref{cor:puente-funcional-amann-sobolev-frontera} replaces the
preceding sequence space by
$\mathbb H^{s,p}_{\partial}(\mathfrak K)$, whose components are modeled on
$\mathbb R^m$ for interior charts and on $\mathbb R^m_+$ for boundary
charts. The operators
$\mathcal S_{\partial,\mathfrak K}$ and
$\mathcal R_{\partial,\mathfrak K}$ satisfy
$\mathcal R_{\partial,\mathfrak K}
 \mathcal S_{\partial,\mathfrak K}=\operatorname{id}$ and estimates
\eqref{eq:comparacion-amann-frontera-analisis}--
\eqref{eq:comparacion-amann-frontera-sintesis}, with positive constants
$c^\partial_{s,p},C^\partial_{s,p},D^\partial_{s,p}$. These depend
only on $m,s,p$, overlap bounds, and uniform bounds on
coordinate changes and cutoff functions through the required order; they do not depend
on the chart or the localized element.
This transfers Euclidean embeddings, multipliers, and local operators
to the charts with uniform constants. This is the role of
Amann's uniform regularity in this book: to provide a coordinate formulation
of bounded geometry adapted to function spaces, rather than an alternative geometry
or a uniformity on a configuration space.

\section{Uniform spaces and Eichhorn uniformities}
\label{sec:espacios-uniformes-eichhorn}

\subsection{Uniform structures}

If $X$ is a set and $U,V\subseteq X\times X$, write
\[
\begin{aligned}
U^{-1}&:=\{(y,x)\mid (x,y)\in U\},\\
U\circ V&:=\{(x,z)\in X\times X\mid
\text{there exists }y\in X\text{ with }(x,y)\in V,\ (y,z)\in U\}.
\end{aligned}
\]
The diagonal of $X\times X$ is denoted by
$\Delta_X:=\{(x,x)\mid x\in X\}$.

\begin{definition}[Uniform structure]
\label{def:estructura-uniforme}\glsadd{estructura-uniforme}
\index{uniform structure}
\index{entourage@\emph{entourage}}
A \textit{uniform structure} on $X$ is a nonempty family
$\varnothing\neq\mathscr U\subseteq\mathcal P(X\times X)$ satisfying:
\begin{enumerate}[label=(\roman*)]
\item if $U\in\mathscr U$ and $U\subseteq V\subseteq X\times X$, then
$V\in\mathscr U$;
\item if $U,V\in\mathscr U$, then $U\cap V\in\mathscr U$;
\item $\Delta_X\subseteq U$ for every $U\in\mathscr U$;
\item if $U\in\mathscr U$, then $U^{-1}\in\mathscr U$;
\item for each $U\in\mathscr U$, there exists $V\in\mathscr U$ such that
$V\circ V\subseteq U$.
\end{enumerate}
The elements of $\mathscr U$ are called \textit{uniform neighborhoods}; they are also
known by the French name \textit{entourages}. The pair $(X,\mathscr U)$ is a \textit{uniform space}.
It is separated if $\displaystyle\bigcap_{U\in\mathscr U}U=\Delta_X$.
\end{definition}

\begin{definition}[Base and cofinality]
\label{def:base-uniformidad}
\index{uniform structure!base}
A family $\mathscr B\subseteq\mathscr U$ is a \textit{base} for the
uniformity $\mathscr U$ if, for each $U\in\mathscr U$, there exists
$B\in\mathscr B$ such that $B\subseteq U$. Two bases $\mathscr B_1$ and
$\mathscr B_2$ are \textit{cofinal} if each element of either
contains an element of the other.

A nonempty family $\mathscr B\subseteq\mathcal P(X\times X)$ generates a uniformity
if all its elements contain $\Delta_X$ and the following conditions
hold: for $B_1,B_2\in\mathscr B$, there exists $B_3\in\mathscr B$ with
$B_3\subseteq B_1\cap B_2$; for $B\in\mathscr B$, there exist
$B',B''\in\mathscr B$ such that $(B')^{-1}\subseteq B$ and
$B''\circ B''\subseteq B$. The generated uniformity consists of the
subsets of $X\times X$ containing some element of $\mathscr B$.
\end{definition}

\begin{definition}[Metrizable uniformity]
\label{def:uniformidad-metrizable}
A uniformity is \textit{metrizable} if a finite-valued metric generates it.
A base is \textit{countable} if it has a countable cofinal
subfamily.
\end{definition}

For $x\in X$ and $U\in\mathscr U$, set
$U[x]:=\{y\in X\mid(x,y)\in U\}$. The sets $U[x]$ form a neighborhood base
for a topology on $X$. The uniform structure contains more
information than this topology: it allows points to be compared without first choosing
a basepoint.

\begin{definition}[Uniform continuity]
\label{def:continuidad-uniforme}
\index{uniform continuity}
Let $(X,\mathscr U_X)$ and $(Y,\mathscr U_Y)$ be uniform spaces. A map
$F\colon X\longrightarrow Y$ is \textit{uniformly continuous} if, for each
$V\in\mathscr U_Y$, there exists $U\in\mathscr U_X$ such that
\[
(F\times F)(U)\subseteq V.
\]
Two uniform structures are equivalent if the identity is uniformly
continuous in both directions.
\end{definition}

\begin{definition}[Filters and uniform completeness]
\label{def:cauchy-completitud-uniforme}
\index{Cauchy filter}
\index{uniform space!complete}
A \textit{filter} on $X$ is a nonempty family $\mathscr F$ of
subsets of $X$ that does not contain $\varnothing$, is closed under
finite intersections, and contains every superset of each of its
elements.
A filter $\mathscr F$ on $X$ is \textit{Cauchy} if, for every
$U\in\mathscr U$, there exists $A\in\mathscr F$ such that $A\times A\subseteq U$.
It converges to $x\in X$ if $U[x]\in\mathscr F$ for every $U\in\mathscr U$.
The uniform space is \textit{complete} if every Cauchy filter converges.
\end{definition}

\begin{lemma}[Filters and sequences for a countable base]
\label{lem:filtros-sucesiones-uniformidad-numerable}
Suppose $\mathscr U$ has a countable base. A sequence $(x_n)$
is called \textit{Cauchy} if, for each $U\in\mathscr U$, there exists $n_0$ such
that $(x_m,x_n)\in U$ whenever $m,n\geq n_0$. Then the space is complete
in the sense of Definition~\ref{def:cauchy-completitud-uniforme} if and
only if every Cauchy sequence converges.
\end{lemma}

\begin{proof}
If every Cauchy filter converges, apply the assumption to the filter of
tails of a Cauchy sequence; convergence of this filter is equivalent to
convergence of the sequence.

Conversely, let $\mathscr F$ be a Cauchy filter and let
$(\mathscr B_k)_{k\in\mathbb N}$ be a countable cofinal base. For each $k$,
there exists $A_k\in\mathscr F$ such that
$A_k\times A_k\subseteq\mathscr B_k$. The sets
$C_k:=A_1\cap\cdots\cap A_k$ belong to $\mathscr F$ and are nonempty;
choose $x_k\in C_k$. Given $U\in\mathscr U$, there exists $j$ with
$\mathscr B_j\subseteq U$. If $m,n\geq j$, then
$x_m,x_n\in C_j\subseteq A_j$ and therefore $(x_m,x_n)\in U$. The sequence
$(x_k)$ is Cauchy and converges, by assumption, to some $x\in X$.

Fix $U\in\mathscr U$. There exists a symmetric entourage $V$ such that
$V\circ V\subseteq U$: intersect an entourage whose square is contained in
$U$ with its inverse. Choose $j$ with $\mathscr B_j\subseteq V$ and then $n\geq j$
with $(x,x_n)\in V$. For $y\in A_j$, we have
$(x_n,y)\in\mathscr B_j\subseteq V$, so $(x,y)\in U$. Consequently,
$A_j\subseteq U[x]$ and $U[x]\in\mathscr F$. This holds for each $U$, so
$\mathscr F$ converges to $x$.
\end{proof}

\begin{proposition}[Metric uniformity]
\label{prop:uniformidad-metrica}
Let $(X,d)$ be a metric space. The sets
\[
U_\varepsilon:=\{(x,y)\in X\times X\mid d(x,y)<\varepsilon\},
\qquad \varepsilon>0,
\]
form a base for a separated uniform structure. Its topology is the
metric topology, its Cauchy sequences are the usual ones, and the uniform
space is complete if and only if $(X,d)$ is complete.
\end{proposition}

\begin{proof}
The diagonal is contained in each $U_\varepsilon$ and
$U_\varepsilon^{-1}=U_\varepsilon$. The triangle inequality gives
$U_{\frac{\varepsilon}{2}}\circ U_{\frac{\varepsilon}{2}}\subseteq U_\varepsilon$.
Moreover, $U_\varepsilon[x]=B_d(x,\varepsilon)$, so the induced topology
is the metric topology. Finally,
$\{x_j\}$ is Cauchy for the uniformity if and only if, given
$\varepsilon>0$, there exists $j_0$ such that
$(x_i,x_j)\in U_\varepsilon$ for $i,j\geq j_0$, which is precisely the
usual metric condition.
Lemma~\ref{lem:filtros-sucesiones-uniformidad-numerable}, applied to the
base $\left(U_{\frac{1}{k}}\right)_{k\geq1}$, identifies completeness by filters with
metric completeness.
\end{proof}

\begin{proposition}[Completion in the metrizable case]
\label{prop:completacion-uniforme-metrizable}
Every separated metrizable uniform space $(X,\mathscr U)$ admits a
metric completion: there exist a separated, metrizable, complete uniform
space $\widehat X$ and a uniform embedding
$\iota\colon X\longrightarrow\widehat X$ with dense image. Here
\textit{uniform embedding} means that $\iota$ is injective and both
$\iota$ and its inverse on $\iota(X)$ are uniformly continuous. The
completion is unique among separated metrizable completions: if
$\jmath\colon X\longrightarrow Y$ is another, there exists a unique uniform
isomorphism $\Phi\colon \widehat X\longrightarrow Y$ such that
$\Phi\circ\iota=\jmath$.
\end{proposition}

\begin{proof}
Choose a metric $d$ generating $\mathscr U$. Consider Cauchy
sequences in $X$ and identify two of them, $(x_j)$ and $(y_j)$, when
$d(x_j,y_j)\to0$. If $[x_j]$ denotes the resulting class, the formula
\[
\widehat d([x_j],[y_j]):=\lim_{j\to\infty}d(x_j,y_j)
\]
is well defined, since
\[
 \bigl|d(x_j,y_j)-d(x_k,y_k)\bigr|
 \leq d(x_j,x_k)+d(y_j,y_k).
\]
This inequality shows that the scalar sequence is Cauchy and its limit
is unchanged when representatives are replaced. Positivity, symmetry, and the
triangle inequality pass to the limit, and
$\widehat d([x_j],[y_j])=0$ is equivalent to the relation defining the classes. Thus $\widehat d$ is a metric. The map
\[
 \iota\colon X\longrightarrow\widehat X,
 \qquad x\longmapsto[(x,x,\ldots)],
\]
is isometric and hence a uniform embedding. If $(x_j)$ represents
$\xi\in\widehat X$, then
\[
 \widehat d(\iota(x_j),\xi)
 =\lim_{k\to\infty}d(x_j,x_k)\longrightarrow0
 \quad\text{as }j\longrightarrow\infty,
\]
so $\iota(X)$ is dense.

Let $(\xi_k)$ be a Cauchy sequence in $\widehat X$. Choose a
subsequence $(\eta_j)$ such that
$\widehat d(\eta_j,\eta_{j+1})<2^{-j-2}$ and, by density, points
$z_j\in X$ such that
$\widehat d(\iota(z_j),\eta_j)<2^{-j-2}$. The triangle inequality gives
$d(z_j,z_{j+1})<2^{-j}$; consequently, for $q>p$,
\[
 d(z_p,z_q)<\sum_{j=p}^{q-1}2^{-j},
\]
and $(z_j)$ is Cauchy. Let $\zeta:=[z_j]$. Then
\[
 \widehat d(\eta_j,\zeta)
 \leq 2^{-j-2}+\sum_{r=j}^{\infty}2^{-r}
 \longrightarrow0.
\]
The subsequence $(\eta_j)$ converges to $\zeta$ and, since $(\xi_k)$ is
Cauchy, we also have $\xi_k\to\zeta$. Thus $\widehat X$ is complete.

For uniqueness, we use the following construction, which also explains
extension to a completion. If $F\colon X\longrightarrow Y$ is uniformly
continuous and $(Y,e)$ is a complete metric space, define, for $\xi=[x_j]$,
\[
 \widehat F(\xi):=\lim_{j\to\infty}F(x_j).
\]
Uniform continuity makes $(F(x_j))$ a Cauchy sequence. If
$d(x_j,y_j)\to0$, given $\varepsilon>0$, choose $\delta>0$ such that
$d(x,y)<\delta$ implies $e(F(x),F(y))<\frac{\varepsilon}{2}$; for large $j$, we
obtain $e(F(x_j),F(y_j))<\frac{\varepsilon}{2}$. Thus the limit is independent of the
representative. The same criterion shows uniform continuity of
$\widehat F$: if $\widehat d(\xi,\eta)<\frac{\delta}{2}$, choose representatives
$(x_j)$ and $(y_j)$ and, for sufficiently large $j$,
$d(x_j,y_j)<\delta$; passage to the limit gives
$e(\widehat F(\xi),\widehat F(\eta))\leq\frac{\varepsilon}{2}<\varepsilon$. The extension is unique
because $Y$ is separated and $\iota(X)$ is dense.

Now let $\jmath\colon X\longrightarrow Y$ be another separated metrizable
completion. The uniform identity between $\iota(X)$ and $\jmath(X)$ extends
to uniformly continuous maps
$\Phi\colon \widehat X\longrightarrow Y$ and
$\Psi\colon Y\longrightarrow\widehat X$. Their compositions fix the
corresponding dense subsets; the preceding uniqueness gives
$\Psi\circ\Phi=\operatorname{id}_{\widehat X}$ and
$\Phi\circ\Psi=\operatorname{id}_Y$. Moreover,
$\Phi\circ\iota=\jmath$, and the same density shows that no other
isomorphism has this property.
\end{proof}

\begin{definition}[Extended metric]
\label{def:metrica-extendida-componente-finita}\glsadd{metrica-extendida}
\index{extended metric}
An \textit{extended metric} on a set $X$ is a map
\[
d\colon X\times X\longrightarrow[0,+\infty]
\]
satisfying the metric axioms, with the arithmetic convention
$a+(+\infty)=+\infty$.
\end{definition}

\begin{definition}[Finite-distance component]
\label{def:componente-uniforme-generalizada}
\index{finite-distance component}
If $d$ is an extended metric and $x\in X$, the \textit{finite-distance
component} of $x$ is
\[
X_x:=\{y\in X\mid d(x,y)<+\infty\}.
\]
The restriction $d\restriction_{X_x}$ is an ordinary metric. The relation
defined by the condition $d(x,y)<+\infty$ is an equivalence relation, and its classes are
the finite-distance components.
\end{definition}

\begin{proposition}[Decomposition into finite-distance components]
\label{prop:componentes-uniformes-equivalencia}
Finite-distance components are open and closed in the topology
generated by balls of finite radius. The completion of $(X,d)$ is the disjoint
union of the metric completions of these components.
\end{proposition}

\begin{proof}
If $y\in X_x$, every ball of finite radius centered at $y$ is contained in
$X_x$ by the triangle inequality; hence $X_x$ is open. Its
complement is the union of the other components and is also open. A
Cauchy sequence eventually lies in a single component, since
for some finite radius all sufficiently late pairs have distance
less than that radius. Consequently, completing each component and
taking the disjoint union produces, and characterizes, the completion of the extended
metric.
\end{proof}

\begin{proposition}[Metrization by chains]
\label{prop:metrizacion-cadenas-uniformidad}
Let $(V_n)_{n\in\mathbb N_0}$ be a decreasing cofinal base of symmetric entourages of
a separated uniform space $X$, and suppose
\[
V_{n+1}\circ V_{n+1}\circ V_{n+1}\subseteq V_n,
\qquad n\geq0.
\]
Define the cost of an edge by
\[
\omega(x,y):=\inf\bigl(\{2^{-n}\mid n\in\mathbb N_0,\ (x,y)\in V_n\}
\cup\{+\infty\}\bigr)
\]
and, for chains $(x_i)_{i=0}^N$ with all vertices in $X$, define the chain
metric by
\begin{equation}
\label{eq:metrizacion-cadenas-uniformidad-abstracta}
d_{\mathrm{cad}}(x,y)
:=\inf\left\{
\sum_{i=0}^{N-1}\omega(x_i,x_{i+1})
\middle|
N\in\mathbb N,\ x=x_0,\ x_N=y
\right\}.
\end{equation}
If there is no chain of finite cost, the infimum is $+\infty$. Then,
for $n\geq1$,
\[
V_n\subseteq\{(x,y)\mid d_{\mathrm{cad}}(x,y)\leq2^{-n}\}
\quad\text{and}\quad
\{(x,y)\mid d_{\mathrm{cad}}(x,y)<2^{-n}\}\subseteq V_{n-1}.
\]
In particular, $d_{\mathrm{cad}}$ is an extended metric generating the
uniformity on each finite-distance component, and
$\min\{1,d_{\mathrm{cad}}\}$ is a finite-valued metric generating the uniformity
on the entire space.
\end{proposition}

\begin{proof}
The first inclusion follows by using a chain with one edge. For
the second, we prove simultaneously for every $n\in\mathbb N$ the following
implication: a chain of $N$ edges with cost less than $2^{-n}$ has
its endpoints in $V_{n-1}$. If $N=1$, the definition of the cost and the
dyadic nature of its values give $(x_0,x_1)\in V_{n+1}\subseteq V_{n-1}$.
An edge of zero cost joins equal points, since the pair belongs to
all $V_j$ and the uniformity is separated.

Now suppose the implication holds for every chain with fewer than $N$ edges.
If the total cost $c$ is zero, all vertices agree. If $c>0$,
choose the first edge for which the accumulated cost reaches or
exceeds $c/2$. The preceding subchain has cost less than $c/2$, and
the following one has cost at most $c/2$. Since $c<2^{-n}$, both costs
are less than $2^{-(n+1)}$. Each subchain has fewer than $N$ edges,
so the induction hypothesis, applied at index $n+1$, places its
endpoints in $V_n$. If a subchain is empty, its endpoints agree and
belong to the diagonal of $V_n$. The central edge has cost less than
$2^{-n}$, so it belongs to $V_{n+1}\subseteq V_n$. Inclusion
$V_n\circ V_n\circ V_n\subseteq V_{n-1}$ completes the induction step.

Symmetry and the triangle inequality for $d_{\mathrm{cad}}$ follow
by reversing and concatenating chains, respectively. If its value is zero, the
second inclusion places the pair in every $V_n$; separation implies that the
points agree. The two inclusions show that balls and
entourages are cofinal. Finally, truncation at one preserves all balls
of radius less than one and assigns distance one to pairs in distinct
components.
\end{proof}

In particular, every separated uniformity with a countable base is metrizable:
intersect each base element with its inverse and refine
recursively to obtain the triple-composition condition in the proposition.

\subsection{Eichhorn uniformities}
\label{subsec:uniformidades-eichhorn}

The preceding uniform regularity is a property of the base manifold and an
atlas. The uniformities introduced now, by contrast, live on
a set of sections, maps, metrics, or diffeomorphisms. Their role is to
compare two objects globally, even when the base is noncompact. The
general uniform framework and its use on complete noncompact manifolds without
boundary are developed in \cite{EichhornUniform2001}; the specific analytic constructions are
formulated below.

\subsubsection{The linear model}

Let $(M,\mathbf g)$ be a complete Riemannian manifold without boundary,
and let $\mathbf E\to M$ be a smooth real vector bundle of finite rank,
with bundle metric $\mathbf h_{\mathbf E}$ and compatible connection
$\nabla:=\nabla^{\mathbf E}$. We use the spaces and norms of
Definition~\ref{Sobolev global y local}, with the connection notation
introduced there. For $1\leq p<\infty$ and $r\in\mathbb N_0$, abbreviate
\[
 \|\mathbf u\|_{p,r;\nabla}
 :=\|\mathbf u\|_{W^{r,p}(M,\mathbf E;\nabla)}.
\]
At integer orders in this chapter,
$H^r(M,\mathbf E;\nabla):=W^{r,2}(M,\mathbf E;\nabla)$.
For tensor bundles, the reference $\mathbf g$ after the semicolon
specifies the tensor metric and connection induced by it.
Derivatives in Sobolev norms are understood weakly; for smooth sections
they agree with classical covariant derivatives.

On the space $L^0(M,\mathbf{E})$ of classes of measurable sections modulo almost everywhere equality,
define the extended metric
\begin{equation}
\label{eq:metrica-extendida-secciones-eichhorn}
 d_{p,r;\nabla}(\mathbf{u},\mathbf{v}):=
 \begin{cases}
  \|\mathbf{u}-\mathbf{v}\|_{p,r;\nabla},
  &\mathbf{u}-\mathbf{v}\in W^{r,p}(M,\mathbf{E}),\\
  +\infty,&\mathbf{u}-\mathbf{v}\notin W^{r,p}(M,\mathbf{E}).
 \end{cases}
\end{equation}

\begin{proposition}[Linear components]
\label{prop:componentes-lineales-eichhorn}
Let $(M,\mathbf{g})$ be a complete Riemannian manifold without boundary, and
let $\mathbf{E}\to M$ be a smooth real vector bundle of finite rank, equipped with a bundle metric $\mathbf{h}_{\mathbf{E}}$ and a compatible connection $\nabla:=\nabla^{\mathbf{E}}$. Fix $1\leq p<\infty$ and
$r\in\mathbb N_0$. Then the sets
\[
 V_\varepsilon:=\{(\mathbf{u},\mathbf{v})\mid
 d_{p,r;\nabla}(\mathbf{u},\mathbf{v})<\varepsilon\},
 \qquad \varepsilon>0,
\]
form a base for a separated metrizable uniformity. The finite-distance component of $\mathbf{u}_0$ is the affine space
$\mathbf{u}_0+W^{r,p}(M,\mathbf{E})$, which is complete when $W^{r,p}(M,\mathbf{E})$ is. If $\nabla'$ is another connection, if the spaces of sections with integrable weak
derivatives of order at most $r$ are the same for both connections, and
if there exist $c_\nabla,C_\nabla>0$ such that
\[
c_\nabla\|\mathbf{u}\|_{p,r;\nabla}
\leq\|\mathbf{u}\|_{p,r;\nabla'}
\leq C_\nabla\|\mathbf{u}\|_{p,r;\nabla}
\qquad (\mathbf{u}\in W^{r,p}(M,\mathbf{E})),
\]
then both connections induce the
same uniformity.
\end{proposition}

\begin{proof}
Symmetry and the triangle inequality for
$d_{p,r;\nabla}$ are inherited from the norm on each finite-distance
component. Hence
$V_{\frac{\varepsilon}{2}}\circ V_{\frac{\varepsilon}{2}}\subseteq V_\varepsilon$;
moreover, $V_\varepsilon^{-1}=V_\varepsilon$,
$\Delta_{L^0(M,\mathbf{E})}\subseteq V_\varepsilon$, and
$V_{\min\{\varepsilon,\delta\}}\subseteq
V_\varepsilon\cap V_\delta$. These are all the conditions in
Definition~\ref{def:base-uniformidad}. The function
$\min\{1,d_{p,r;\nabla}\}$ is a finite-valued metric generating the same
entourages of radius less than one, proving metrizability. On the other
hand,
$d_{p,r;\nabla}(\mathbf{u}_0,\mathbf{u})<\infty$ is equivalent to
$\mathbf{u}-\mathbf{u}_0\in W^{r,p}(M,\mathbf{E})$, identifying the component. The translation
$\mathbf{u}\mapsto \mathbf{u}-\mathbf{u}_0$ is an isometry from this component onto $W^{r,p}(M,\mathbf{E})$,
which gives completeness. Finally, two norms satisfying the
explicit inequalities in the statement have
the same uniform balls after changing the radius by the factors
$c_\nabla,C_\nabla$ and hence cofinal bases of entourages.
\end{proof}

\begin{example}[Connections cannot be changed arbitrarily]
\label{ej:conexion-no-acotada-sobolev}
On the trivial bundle $\mathbb R\times\mathbb R^2$ with its Euclidean metric,
consider the metric connections $\nabla=d$ and $\nabla'=d+xJ\,dx$, where
$Je_1=e_2$ and $Je_2=-e_1$.
The field $\mathbf u(x)=(1+x^2)^{-1/2}e_1$ belongs to $H^1(\mathbb R,\mathbb R^2)$ for $\nabla$:
both $(1+x^2)^{-1/2}$ and its derivative $-x(1+x^2)^{-3/2}$ belong to $L^2(\mathbb R)$.
For every $j\in\mathbb N_0$, successive differentiation gives
\[
 \frac{d^j}{dx^j}(1+x^2)^{-1/2}
 =\frac{P_j(x)}{(1+x^2)^{j+1/2}},\qquad
 P_0=1,\quad P_{j+1}=(1+x^2)P'_j-(2j+1)xP_j.
\]
The recurrence implies $\deg P_j\leq j$. If $C_j$ is the sum of the
absolute values of its coefficients, then
\[
 \left|\frac{d^j}{dx^j}(1+x^2)^{-1/2}\right|
 \leq C_j|x|^{-j-1},\qquad |x|\geq1.
\]
The square of this bound is integrable on $\{x\in\mathbb R\mid |x|\geq1\}$;
on $[-1,1]$ the derivative is continuous. Thus $\mathbf u$ belongs to
$H^r(\mathbb R,\mathbb R^2)$ defined by $\nabla$ for every integer $r\geq1$.
However,
\[
 \nabla'\mathbf u=
 \bigl(-x(1+x^2)^{-3/2}e_1+x(1+x^2)^{-1/2}e_2\bigr)\,dx
\]
does not belong to $L^2(\mathbb R,\mathbb R^2)$, since the square of its second component tends to one.
Both connections are compatible with the same metric and even have zero
curvature, because the base is one-dimensional. What fails is uniform control
of their difference. Thus equivalence of spaces and norms in the preceding
proposition is a hypothesis, not a consequence of smoothness or bounded
curvature of the two connections.
\end{example}

\subsubsection{Effective criteria for changing the connection}
\label{subsubsec:comparacion-conexiones-eichhorn}

The equivalence hypothesis in
Proposition~\ref{prop:componentes-lineales-eichhorn} can be checked
using the difference of the connections. We distinguish two criteria:
uniform bounds on its derivatives and its Sobolev integrability.
The latter is the setting of Theorem~3.22 in Chapter I of
\cite{Eichhorn2007}. We develop the estimates, which also specify
which reference connections can be interchanged.

\begin{lemma}[The module property at integer orders]
\label{lem:modulo-sobolev-uniforme-entero}
Let $M$ be a Riemannian manifold without boundary, of dimension
$m\geq1$, and let $\mathbf E_1,\mathbf E_2,\mathbf E_3$ be finite-rank metric bundles.
Let $*\colon\mathbf E_1\otimes\mathbf E_2\to\mathbf E_3$ be the tensor product
or a parallel contraction with uniformly bounded fiber norm.
Assume that the base and the bundles have bounded geometry of order
$k\geq a$. If
$1<p<\infty$, $a,b\in\mathbb N_0$, $a>m/p$ and $0\leq b\leq a$, then
the tensor product, possibly followed by a metric contraction,
satisfies
\begin{equation}\label{eq:modulo-sobolev-uniforme-entero}
 \|\mathbf u*\mathbf v\|_{W^{b,p}(M,\mathbf E_3)}
 \leq C\|\mathbf u\|_{W^{a,p}(M,\mathbf E_1)}\|\mathbf v\|_{W^{b,p}(M,\mathbf E_2)}.
\end{equation}
The constant depends on the orders, the ranks and the geometric bounds,
not on the sections.
\end{lemma}

\begin{proof}
We use the trivializations and the equivalence with covariant norms
from Proposition~\ref{prop:covariantes-localizacion-orden-entero}.
Recall how local estimates give global ones. For
$q\geq p$, a local embedding estimate whose constant is independent of
the chart is summed using
$\|(c_\alpha)\|_{\ell^q(I)}\leq\|(c_\alpha)\|_{\ell^p(I)}$; for $q=\infty$
we take the supremum. The bounded multiplicity of the covering controls
the sums of the $p$th powers of the local norms. Thus the Euclidean
embeddings already proved give, for the tensor bundles occurring here,
\[
 W^{d,p}(M,\mathbf F)\hookrightarrow L^q(M,\mathbf F)\quad\text{if }q\geq p,
 \quad d>m\left(\frac1p-\frac1q\right),
 \qquad W^{d,p}(M,\mathbf F)\hookrightarrow L^\infty(M,\mathbf F)\quad\text{if }d>m/p.
\]
Only orders $d\leq a$ are used, so bounded geometry of all orders is
not required.

By Leibniz it suffices to estimate, for $0\leq i\leq j\leq b$,
$\nabla^i\mathbf u*\nabla^{j-i}\mathbf v$ in $L^p(M,(T^*M)^{\otimes j}\otimes \mathbf E_3)$.
If $i=0$, use $\|\mathbf u\|_\infty\leq C\|\mathbf u\|_{W^{a,p}(M,\mathbf E_1)}$.
If $a-i=0$, necessarily $i=j=b=a$; in this case use
$\|\mathbf v\|_\infty\leq C\|\mathbf v\|_{W^{a,p}(M,\mathbf E_2)}$.
In every other case, the integers $d_1=a-i$ and $d_2=b-j+i$ are positive
and $d_1+d_2=a+b-j\geq a>m/p$. Consequently,
\[
 A:=\max\left\{0,\frac1p-\frac{d_1}{m}\right\},\qquad
 B:=\max\left\{0,\frac1p-\frac{d_2}{m}\right\}
 \quad\text{satisfy}\quad A+B<\frac1p.
\]
We can choose $t_1,t_2\in(0,1/p)$ with $t_1>A$, $t_2>B$ and
$t_1+t_2=1/p$: choose $t_1$ in the nonempty interval
$(A,1/p-B)$, avoiding its endpoints. Setting $p_\nu=1/t_\nu$, the
preceding embeddings and Hölder give
\[\begin{aligned}
 \|\nabla^i\mathbf u*\nabla^{j-i}\mathbf v\|_{L^p(M,(T^*M)^{\otimes j}\otimes \mathbf E_3)}\\
 &\leq C\|\nabla^i\mathbf u\|_{L^{p_1}(M,(T^*M)^{\otimes i}\otimes \mathbf E_1)}
          \|\nabla^{j-i}\mathbf v\|_{L^{p_2}(M,(T^*M)^{\otimes(j-i)}\otimes \mathbf E_2)}\\
 &\leq C\|\mathbf u\|_{W^{a,p}(M,\mathbf E_1)}\|\mathbf v\|_{W^{b,p}(M,\mathbf E_2)}.
\end{aligned}\]
We sum the finitely many Leibniz terms. Local approximation by smooth
sections, together with the uniform estimates, gives the same identity
and bound for weak derivatives. Metric contractions are parallel and
bounded in each fiber, and do not change the argument.
\end{proof}

\begin{proposition}[A difference with uniformly bounded derivatives]
\label{prop:cambio-conexion-cotas-uniformes}
On a fixed Riemannian base, let $\nabla$ and $\nabla'$ be smooth
connections on the same metric bundle $\mathbf E$, and write
$\mathbf a=\nabla'-\nabla$. Let $r\geq1$ be an integer. If
\[
 \displaystyle\max_{0\leq j\leq r-1}
               \|\nabla^j\mathbf a\|_{L^\infty(M,(T^*M)^{\otimes j+1}\otimes \operatorname{End}\mathbf E)}<\infty,
\]
then, for $1\leq p<\infty$ and every integer $0\leq s\leq r$, the spaces
$W^{s,p}(M,\mathbf E;\nabla)$ and $W^{s,p}(M,\mathbf E;\nabla')$
coincide with equivalent norms. No separate assumption that the
$\nabla'$-derivatives of $-\mathbf a$ are bounded is needed.
\end{proposition}

\begin{proof}
On the bundles $T^{(0,j)}M\otimes\mathbf E$ we keep the same
Levi--Civita connection on the base factors. The induced difference acts
only on the $\mathbf E$ factor. Iterating
$\nabla'=\nabla+\mathbf a*$ and applying Leibniz expresses every term of
$(\nabla')^s\mathbf u-\nabla^s\mathbf u$ as a contraction of
\[
 (\nabla^{i_1}\mathbf a)\otimes\cdots\otimes
 (\nabla^{i_\ell}\mathbf a)\otimes\nabla^j\mathbf u,
 \qquad \ell\geq1,\quad
 j+\displaystyle\sum_{\nu=1}^{\ell}(i_\nu+1)=s.
\]
This description follows by induction: $\nabla$ increases the order of
one factor by one, whereas $\mathbf a*$ adds a factor with $i_\nu=0$.
In particular, $i_\nu\leq s-1$ and $j\leq s-1$.
The bounds on the coefficients give the continuous inclusion from the
unprimed space into the primed space.

On $T^*M\otimes\operatorname{End}\mathbf E$, the induced difference
is the commutator action of $\mathbf a$ on the endomorphism factor.
The same induction expresses $(\nabla')^j\mathbf a$ in terms of products
of $\nabla^{i}\mathbf a$, with $i\leq j$. All these expressions are
bounded for $j\leq r-1$. We now apply the first argument with
$\nabla'$ as reference and difference $-\mathbf a$, obtaining the
opposite continuous inclusion. These are identities of differential
operators with smooth coefficients, so they hold in distributions.
They therefore compare the full domains, not merely the norms on common
smooth sections. The case $s=0$ is immediate.
\end{proof}

\begin{theorem}[A connection difference of Sobolev regularity]
\label{teo:cambio-conexion-diferencia-sobolev}
Let $(M,\mathbf g)$ be complete and without boundary, of dimension $m$,
and let $\mathbf E\to M$ be a finite-rank metric bundle. Let
$\nabla,\nabla'$ be smooth metric connections. Assume the bounded
geometry conditions for the base and for both connections up to an
order $k\geq r$, where
\[
 1<p<\infty,\qquad r\in\mathbb N,\qquad r>\frac mp+1.
\]
If $\mathbf a=\nabla'-\nabla$ belongs to
$W^{r-1,p}(M,T^*M\otimes\operatorname{End}\mathbf E;\nabla)$,
then, for every integer $0\leq s\leq r$,
\begin{equation}\label{eq:estimacion-cambio-conexion-sobolev}
 W^{s,p}(M,\mathbf E;\nabla)=W^{s,p}(M,\mathbf E;\nabla'),
 \qquad
 c_s\|\mathbf u\|_{p,s;\nabla}
 \leq\|\mathbf u\|_{p,s;\nabla'}
 \leq C_s\|\mathbf u\|_{p,s;\nabla}.
\end{equation}
Moreover, under these same geometric assumptions,
\begin{equation}\label{eq:diferencia-conexion-simetrica}
 \mathbf a\in W^{r-1,p}(M,T^*M\otimes\operatorname{End}\mathbf E;\nabla)
 \quad\Longleftrightarrow\quad
 \mathbf a\in W^{r-1,p}(M,T^*M\otimes\operatorname{End}\mathbf E;\nabla').
\end{equation}
\end{theorem}

\begin{proof}
Put $a=r-1>m/p$ and $A=\|\mathbf a\|_{p,a;\nabla}$.
The module lemma shows that multiplication by $\mathbf a$ is a bounded
operator on $W^{b,p}(M,\mathbf E;\nabla)$ for every integer $0\leq b\leq a$,
with norm at most $C_b A$. The same holds on each tensor bundle of
order at most $r$ occurring in the derivatives.

Let $\mathbf u\in W^{s,p}(M,\mathbf E;\nabla)$, $1\leq s\leq r$.
We prove by induction, for $0\leq j\leq s$, that
\[
 (\nabla')^j\mathbf u\in W^{s-j,p}(M,(T^*M)^{\otimes j}\otimes \mathbf E;\nabla),\qquad
 \|(\nabla')^j\mathbf u\|_{p,s-j;\nabla}
 \leq C_j(1+A)^j\|\mathbf u\|_{p,s;\nabla}.
\]
There is nothing to prove for $j=0$. If the assertion holds for
$j<s$, write $\mathbf v=(\nabla')^j\mathbf u$. Then
$\nabla'\mathbf v=\nabla\mathbf v+\mathbf a*\mathbf v$.
The first term belongs to $W^{s-j-1,p}(M,(T^*M)^{\otimes j+1}\otimes \mathbf E;\nabla)$ by definition;
the second does as well, by the lemma, since
$0\leq s-j-1\leq r-1=a$. This proves the induction step.
In particular, all primed derivatives up to order $s$ are in $L^p(M,(T^*M)^{\otimes j}\otimes \mathbf E)$,
giving the upper bound in
\eqref{eq:estimacion-cambio-conexion-sobolev}.

We must justify the opposite bound without presupposing the desired
conclusion. Apply only the inclusion just established to the bundle
$T^*M\otimes\operatorname{End}\mathbf E$, to the section $\mathbf a$
and to the order $s=a$. The difference of its induced connections is
a fixed sum of left and right actions of $\mathbf a$; its multiplication
norm is therefore at most $C A$. The same argument yields
\[
 \|\mathbf a\|_{p,a;\nabla'}
 \leq C(1+A)^a A<\infty.
\]
This proves one implication of
\eqref{eq:diferencia-conexion-simetrica}. Since the geometric hypotheses
were imposed on both connections, the module lemma also applies with
$\nabla'$ as reference. Repeat the direct inclusion with reference
$\nabla'$ and difference $-\mathbf a$. This gives the reverse inclusion
of the spaces of order $s$, its continuous estimate, and the other
implication in \eqref{eq:diferencia-conexion-simetrica}.

All preceding operations make sense for weak derivatives: $\nabla$
acts in distributions and the products lie in the stated spaces by the
module lemma. The proof thus applies to any section of either Sobolev
domain. It does not use a prior identification of those domains.
\end{proof}

\begin{corollary}[Connection perturbations with Sobolev coefficients]
\label{cor:conexion-sobolev-dominios-debiles}
Let $1<p<\infty$ and let $r>m/p+1$ be an integer. Suppose that the base
and the reference connection $\nabla$ admit the uniform localizations
and the module property of Lemma~\ref{lem:modulo-sobolev-uniforme-entero}
through order $r-1$, in the tensor bundles used. It suffices
that $\mathbf a\in W^{r-1,p}(M,T^*M\otimes\operatorname{End}\mathbf E;\nabla)$,
without assuming smoothness of $\mathbf a$ or bounded geometry of
$\nabla'=\nabla+\mathbf a$, to compare the following weak domains.
For a section $\mathbf u\in L^p(M,\mathbf E)$ define recursively
\[
 \mathbf v_0=\mathbf u,\qquad
 \mathbf v_{j+1}=\nabla\mathbf v_j+\mathbf a*\mathbf v_j.
\]
As long as $\mathbf v_j\in L^p(M,(T^*M)^{\otimes j}\otimes \mathbf E)$, the right-hand side is a well-defined
distribution, since $\mathbf a$ is bounded. The domain on which
$\mathbf v_j\in L^p(M,(T^*M)^{\otimes j}\otimes\mathbf E)$ for
every $j\in\{0,\ldots,s\}$, with $0\leq s\leq r$,
is exactly $W^{s,p}(M,\mathbf E;\nabla)$, and the norm
$\displaystyle\sum_{j=0}^s\|\mathbf v_j\|_{L^p(M,(T^*M)^{\otimes j}\otimes \mathbf E)}$ is equivalent to its usual norm.
\end{corollary}

\begin{proof}
The inclusion of $W^{s,p}(M,\mathbf E;\nabla)$ in the stated domain and its
estimate were proved by the direct induction in the preceding theorem;
only weak derivatives and the module lemma for $\nabla$ are used.
For the reverse inclusion, assume $\mathbf v_j\in L^p(M,(T^*M)^{\otimes j}\otimes \mathbf E)$ for
$0\leq j\leq s$. We prove by induction on $\ell$ that
\[
 \mathbf v_j\in W^{\ell,p}(M,(T^*M)^{\otimes j}\otimes \mathbf E;\nabla)
 \quad(0\leq j\leq s-\ell),\qquad 0\leq\ell\leq s.
\]
The case $\ell=0$ is the hypothesis. If the assertion holds for
$\ell-1$, then, for $j\leq s-\ell$,
\[
 \mathbf v_j\in W^{\ell-1,p}(M,(T^*M)^{\otimes j}\otimes\mathbf E;\nabla),
 \qquad
 \mathbf v_{j+1}\in W^{\ell-1,p}(M,(T^*M)^{\otimes(j+1)}\otimes\mathbf E;\nabla).
\] By the module lemma,
$\mathbf a*\mathbf v_j$ belongs to the second space, since
$\ell-1\leq r-1$. The distributional identity
\[
 \nabla\mathbf v_j=\mathbf v_{j+1}-\mathbf a*\mathbf v_j
\]
implies $\mathbf v_j\in W^{\ell,p}(M,(T^*M)^{\otimes j}\otimes \mathbf E;\nabla)$. Moreover,
\[
 \|\mathbf v_j\|_{p,\ell;\nabla}
 \leq C\left(\|\mathbf v_j\|_{L^p(M,(T^*M)^{\otimes j}\otimes \mathbf E)}
       +\|\mathbf v_{j+1}\|_{p,\ell-1;\nabla}
       +\|\mathbf a\|_{p,r-1;\nabla}
                         \|\mathbf v_j\|_{p,\ell-1;\nabla}\right).
\]
The induction also bounds each term by a constant times
$\displaystyle\sum_{i=0}^{s}\|\mathbf v_i\|_{L^p(M,(T^*M)^{\otimes i}\otimes \mathbf E)}$. Taking $j=0$ and $\ell=s$
gives the reverse inclusion and estimate.
\end{proof}

Consequently, connections satisfying either criterion induce the same
linear uniformity. The expression ``admissible connections'' must be
understood with this control of their difference; bounded geometry of
each connection alone is not a substitute, as shown by
Example~\ref{ej:conexion-no-acotada-sobolev}. For two connections on a
fixed bundle, their difference is a section of
$T^*M\otimes\operatorname{End}\mathbf E$; the extended distance
\eqref{eq:metrica-extendida-secciones-eichhorn} applies to that section
once a reference in this class has been fixed.

\subsubsection{Maps and exponential comparison}

We use conditions $(I)$ and $(B_k)$ from
Definition~\ref{def:geometria-acotada-orden-k}. The constants in the
estimates will depend on the lower injectivity-radius bound $i_0$,
the curvature bounds $K_0,\ldots,K_k$, and the other data specified
in each statement.

In what follows, $(M,\mathbf{g})$ and $(N,\mathbf{h}_N)$ denote complete manifolds without boundary
satisfying the conditions in
Definition~\ref{def:geometria-acotada-orden-k} to the required order.

If $\mathbf{E}\longrightarrow M$ is a smooth real vector bundle of finite rank, equipped with a bundle metric $\mathbf{h}_{\mathbf{E}}$ and a connection $\nabla:=\nabla^{\mathbf{E}}$, and
$j\in\mathbb N_0$, set
\[
C_b^j(M,\mathbf{E})
:=\left\{\mathbf{u}\in C^j(M,\mathbf{E})\middle|
\|\mathbf{u}\|_{C_b^j(M,\mathbf{E})}:=
\sum_{\ell=0}^j\sup_{x\in M}|\nabla^\ell \mathbf{u}(x)|_{\mathbf{g},\mathbf{h}_{\mathbf{E}}}<+\infty
\right\}.
\]
For $j=0$, only continuity and boundedness of $\mathbf{u}$ are required. On the pullback bundle
$f^*TN$, use the connection induced by $f$ and the
Levi--Civita connection of $N$.

For $q\in\mathbb N$, denote by $C_b^{\infty,q}(M,N)$ the set of
smooth maps $f\colon M\longrightarrow N$ such that
$\|\nabla^jdf\|_\infty<+\infty$ for $0\leq j\leq q-1$. When uniform
constants are needed, fix
$\mathbf{B}=(B_0,\ldots,B_{q-1})\in(0,+\infty)^q$ and write
\[
 C_{b,\mathbf{B}}^{\infty,q}(M,N)
 :=\left\{f\in C_b^{\infty,q}(M,N)
 \middle|
 \|\nabla^jdf\|_\infty\leq B_j,\ 0\leq j\leq q-1\right\}.
\]
The constants will be common to all maps in a class with
$\mathbf B$ fixed. Allowing each map its own bound does not give
a common constant; this distinction matters when checking the axioms
of a uniformity.

If $\mathbf{X}\in\Gamma(f^*TN)$ and
$\|\mathbf{X}\|_{C^0(M,f^*TN;f^*\mathbf{h}_N)}<i_0(N)$, set
$(\operatorname{Exp}_f\mathbf{X})(x):=\exp_{f(x)}\mathbf{X}(x)$ and measure its derivatives with the
induced connection on $f^*TN$ by
\[
\|\mathbf{X}\|_{b,r;f}:=\sum_{j=0}^{r}
\|\nabla^j\mathbf X\|_\infty.
\]
For
$0<\varepsilon<i_0(N)/2$, define
\begin{equation}
\label{eq:entourage-mapeos-eichhorn-b}
 \mathcal V^{b,r}_\varepsilon
 :=\left\{(f,h)\middle|
 h=\operatorname{Exp}_f\mathbf{X}\text{ for some }\mathbf{X}\in\Gamma(f^*TN),\
 \|\mathbf{X}\|_{b,r;f}<\varepsilon\right\}.
\end{equation}

Now suppose $1<p<\infty$ and $r\in\mathbb N$ with $r>\dim(M)/p+1$, and fix
$0<\rho<i_0(N)/2$. In the Sobolev regime, we use
\begin{equation}
\label{eq:entourage-mapeos-eichhorn-pr}
 \mathcal V^{p,r}_{\varepsilon,\rho}
 :=\left\{(f,h)\middle|
 \begin{array}{l}
 h=\operatorname{Exp}_f\mathbf{X}\text{ for some }
 \mathbf{X}\in W^{r,p}(M,f^*TN),\\
 \|\mathbf{X}\|_{W^{r,p}(M,f^*TN)}<\varepsilon,
 \quad \|\mathbf{X}\|_{C^0(M,f^*TN;f^*\mathbf{h}_N)}<\rho
 \end{array}\right\}.
\end{equation}
The pointwise condition keeps the comparison inside a common
normal ball. For $q\geq r$, on $C_{b,\mathbf{B}}^{\infty,q}(M,N)$, localization
results and Sobolev embedding give a constant
$C_S=C_S(M,N,p,r,\rho,\mathbf{B})>0$ such that
\[
\|\mathbf{X}\|_{C^0(M,f^*TN;f^*\mathbf{h}_N)}
\leq C_S\|\mathbf{X}\|_{W^{r,p}(M,f^*TN)}.
\]
Thus, in a class with fixed
data, the condition $\|\mathbf{X}\|_{C^0(M,f^*TN;f^*\mathbf{h}_N)}<\rho$ follows from
$\varepsilon<\rho/C_S$.

The central property is that inversion and composition of exponential
comparisons preserve control. If
$h=\operatorname{Exp}_f\mathbf{X}$, denote by
$P_{\mathbf{X}}\colon f^*TN\longrightarrow h^*TN$ parallel transport along
$t\mapsto\exp_{f(x)}(t\mathbf{X}(x))$.

\begin{lemma}[Uniform comparison estimates]
\label{lem:estimaciones-comparacion-eichhorn}
Let $(M,\mathbf{g})$ and $(N,\mathbf{h}_N)$ be complete Riemannian
manifolds without boundary satisfying $(I)$ and $(B_k)$ to the required
order. Fix a radius $\rho>0$ small enough that balls of radius $4\rho$
lie in uniform normal domains and the vertical differentials
of the exponential map and their inverses have common bounds. It may
be chosen in terms of $i_0(N)$ and the curvature bound, with
$4\rho<i_0(N)$. Also fix a profile $\mathbf{B}$ bounding the derivatives of
the maps under consideration through the necessary order. Throughout the
statement, the comparison sections satisfy
\[
\|\mathbf X\|_{C^0(M,f^*TN;f^*\mathbf h_N)}<\rho,
\qquad
\|\mathbf Y\|_{C^0(M,h^*TN;h^*\mathbf h_N)}<\rho
\]
whenever the corresponding section occurs. For each $0\leq j\leq r$,
there exist continuous polynomials $P_j,Q_j$ with nonnegative coefficients and no
constant term, depending only on the dimensions, $j$, $\rho$,
the uniform geometric data of $M$ and $N$ through order $r$,
and $\mathbf{B}$, and satisfying:
\begin{enumerate}[label=(\roman*)]
\item if $h=\operatorname{Exp}_f\mathbf{X}$, then
\[
\begin{aligned}
&\|\nabla^j(P_{\mathbf{X}}\mathbf{X})\|_{
C^0(M,T^{(0,j)}(TM)\otimes h^*TN;\mathbf{g},h^*\mathbf{h}_N)}\\
&\quad\leq P_j\!\left(
\|\mathbf{X}\|_{C^0(M,f^*TN;f^*\mathbf{h}_N)},\ldots,
\|\nabla^j\mathbf{X}\|_{C^0(M,T^{(0,j)}(TM)\otimes f^*TN;\mathbf{g},f^*\mathbf{h}_N)}
\right);
\end{aligned}
\] since
$f=\operatorname{Exp}_h(-P_{\mathbf{X}}\mathbf{X})$, this bound controls inversion;
\item if $h=\operatorname{Exp}_f\mathbf{X}$ and
$\ell=\operatorname{Exp}_h\mathbf{Y}$, there exists a unique
$\mathbf Z\in\Gamma(f^*TN)$ with
$\|\mathbf Z\|_{C^0(M,f^*TN;f^*\mathbf h_N)}<2\rho$ and
$\ell=\operatorname{Exp}_f\mathbf Z$. It is given by the inverse of the
exponential map on the normal ball,
$\mathbf Z(x)=\exp_{f(x)}^{-1}(\ell(x))$, and
$\|\nabla^j\mathbf{Z}\|_{
C^0(M,T^{(0,j)}(TM)\otimes f^*TN;\mathbf{g},f^*\mathbf{h}_N)}
\leq Q_j(\|\mathbf{X}\|_{b,j;f},\|\mathbf{Y}\|_{b,j;h})$.
\end{enumerate}
We may take $P_0(t)=t$ and $Q_0(a,b)=a+b$. Moreover, given $\eta>0$, there exists
$\delta>0$, depending only on the preceding data, such that
$Q_j(a,b)<\eta$ whenever $a+b<\delta$. If
$r>\dim(M)/p+1$, there exist polynomials $\widehat P_r,\widehat Q_r$ with the
same dependencies, together with the Sobolev module and embedding constants,
for which
\[
\begin{aligned}
\|P_{\mathbf{X}}\mathbf{X}\|_{W^{r,p}(M,h^*TN)}
&\leq\widehat P_r(\|\mathbf{X}\|_{W^{r,p}(M,f^*TN)}),\\
\|\mathbf{Z}\|_{W^{r,p}(M,f^*TN)}
&\leq\widehat Q_r\!\left(
\|\mathbf{X}\|_{W^{r,p}(M,f^*TN)},
\|\mathbf{Y}\|_{W^{r,p}(M,h^*TN)}\right).
\end{aligned}
\]
\end{lemma}

\begin{proof}
The pointwise bounds give
\[
\begin{aligned}
d_{\mathbf h_N}(f(x),\ell(x))
&\leq d_{\mathbf h_N}(f(x),h(x))+d_{\mathbf h_N}(h(x),\ell(x))\\
&\leq |\mathbf X(x)|+|\mathbf Y(x)|<2\rho<i_0(N).
\end{aligned}
\]
The inverse of the exponential map on the normal ball therefore defines
$\mathbf Z$ smoothly. Its norm is the distance from $f(x)$ to $\ell(x)$;
taking the supremum proves the stated bound, and injectivity on the
normal ball proves uniqueness in that class.

At order zero, parallel transport is isometric, so
$|P_{\mathbf{X}}\mathbf{X}|_{h^*\mathbf{h}_N}=|\mathbf{X}|_{f^*\mathbf{h}_N}$. Together with the distance estimate, this proves the order-zero bounds
with $P_0(t)=t$ and $Q_0(a,b)=a+b$.

At first order, let $x(s)$ be a curve with $x(0)=x$ and $x'(0)=\mathbf{V}$, and
consider
$c(s,t)=\exp_{f(x(s))}(t\mathbf{X}(x(s)))$. The field
$\mathbf{J}(t)=\partial_sc(0,t)$ satisfies
\begin{equation}
\label{eq:jacobi-comparacion-eichhorn}
 \frac{D^2\mathbf{J}}{dt^2}
 +\operatorname{Rm}_{\mathbf{h}_N}(\mathbf{J},\dot\gamma)\dot\gamma=0,
 \qquad \mathbf{J}(0)=df(\mathbf{V}),\qquad
 \frac{D\mathbf{J}}{dt}(0)=\nabla_{\mathbf{V}}\mathbf{X}.
\end{equation}
In a parallel frame, the integral form of this equation and
Grönwall's lemma~\ref{lema: gronwall} imply
$\displaystyle \sup_{t\in[0,1]}|\mathbf{J}(t)|_{\mathbf{h}_N}\leq
 e^{K_0(N)\rho^2}
 (|df(\mathbf{V})|_{\mathbf{h}_N}+|\nabla_{\mathbf{V}}\mathbf{X}|_{\mathbf{h}_N})$. Integrating the Jacobi equation once
gives
\begin{equation}
\label{eq:cota-derivada-jacobi-eichhorn}
 \sup_{t\in[0,1]}\left|\frac{D\mathbf{J}}{dt}(t)\right|_{\mathbf{h}_N}
 \leq |\nabla_{\mathbf{V}}\mathbf{X}|_{\mathbf{h}_N}
 +K_0(N)|\mathbf{X}|_{f^*\mathbf{h}_N}^2e^{K_0(N)\rho^2}
 \bigl(|df(\mathbf{V})|_{\mathbf{h}_N}+|\nabla_{\mathbf{V}}\mathbf{X}|_{\mathbf{h}_N}\bigr).
\end{equation}
Since the connection is torsion-free,
$\nabla_{\mathbf{V}}(P_{\mathbf{X}}\mathbf{X})=D_t\mathbf{J}(1)$, and the preceding bound proves the case $j=1$.

To compose two exponential comparisons, write
$\widetilde{\mathbf{Y}}:=P_{\mathbf{X}}^{-1}\mathbf{Y}\in\Gamma(f^*TN)$. On the common normal domain,
consider the smooth map
\[
 \Psi(p,v,y):=\exp_p^{-1}\!\left(\exp_{\exp_pv}(P_vy)\right),
 \qquad p\in N,\quad v,y\in T_pN,
\]
where $P_v\colon T_pN\longrightarrow T_{\exp_pv}N$ is parallel
transport along $t\mapsto\exp_p(tv)$. The section directly comparing
$\ell=\operatorname{Exp}_h\mathbf{Y}$ with $f$ is then
$\mathbf{Z}(x)=\Psi(f(x),\mathbf{X}(x),\widetilde{\mathbf{Y}}(x))$. Differentiating this identity
produces the preceding Jacobi fields and the linear equation for parallel
transport. In particular, there exists a constant $A_1$ depending only on
$\dim N$, $\rho$, and $K_0(N)$ such that
\[
 |\nabla \mathbf{Z}|_{\mathbf{g},f^*\mathbf{h}_N}\leq A_1\bigl[
 (1+|\mathbf{X}|_{f^*\mathbf{h}_N}+|\widetilde{\mathbf{Y}}|_{f^*\mathbf{h}_N})
 (|\nabla \mathbf{X}|_{\mathbf{g},f^*\mathbf{h}_N}
  +|\nabla\widetilde{\mathbf{Y}}|_{\mathbf{g},f^*\mathbf{h}_N})
 +( |\mathbf{X}|_{f^*\mathbf{h}_N}+|\widetilde{\mathbf{Y}}|_{f^*\mathbf{h}_N})
 |df|_{\mathbf{g},f^*\mathbf{h}_N}\bigr].
\]
Derivatives of $\widetilde{\mathbf{Y}}=P_{\mathbf{X}}^{-1}\mathbf{Y}$ are estimated using the same
parallel transport equation; their bounds depend polynomially on the
derivatives of $\mathbf{X}$ and $\mathbf{Y}$ and introduce no constant term.

To make the induction step precise, suppose bounds have been obtained through
order $j-1$, with $2\leq j\leq r$. Differentiate equation
\eqref{eq:jacobi-comparacion-eichhorn} a total of $j-1$ times. The recurrence
\eqref{eq:indexacion-exacta-fuentes-jacobi}--
\eqref{eq:peso-exacto-fuentes-jacobi} separates the order-$j$ unknown
from the source: the operator acting on this unknown remains the Jacobi
operator, and each source term contains only variation fields
of lower orders and curvature derivatives of the orders indicated there.
Proposition~\ref{prop: leibniz operador *} allows the
contractions to be differentiated without changing their nature. Differentiating
a curvature factor raises its order by one; differentiating a
field raises the order of that field by one. These are exactly the two
operations in the recurrence, and the weight identity prevents the source
from containing a second unknown of order $j$. Its factors are therefore
controlled by the induction hypothesis and the fixed data.

We must also preserve vanishing when $\mathbf X=0$.
To do so, subtract from the differentiated equation the one corresponding to the
constant variation $c(x,t)=f(x)$. The terms containing only
$df,\ldots,\nabla^{j-1}df$ cancel. Every remaining term contains at
least one factor among $\mathbf X,\nabla\mathbf X,\ldots,\nabla^j\mathbf X$.
Integrating the linear equation and applying Grönwall gives the order-$j$
bound for the final vector $P_{\mathbf X}\mathbf X$, with a polynomial having no
constant term. This completes the step $j-1\longrightarrow j$, starting from the
zero and first orders already computed.

For composition, proceed with the map $\Psi$ and the parallel
transport equation. Here subtract the identity $\Psi(p,0,0)=0$. Differentiating it
with respect to $p$ again eliminates all terms containing only derivatives of the center.
The chain rule of order $j$ leaves products of
derivatives of $\mathbf X$ and $\mathbf Y$ whose total order is at most
$j$; their coefficients are bounded by the same Jacobi recurrence and
the bound for the inverse differential of the exponential map. The sum of
their bounds defines $Q_j$ and retains the absence of a constant term.

To pass to Sobolev norms, use the uniform trivializations and
the norm equivalence of
Proposition~\ref{prop:covariantes-localizacion-orden-entero}. Each term
has at least one integrable factor. Products of derivatives of the
sections are estimated by Hölder and Sobolev embeddings, assigning to
each factor the order it retains after differentiation; the sum of the
differentiation orders does not exceed $r$. The order-zero factors are
controlled by embedding into $C_b^1(M,f^*TN)$, and the coefficient derivatives by
the geometric data and $\mathbf B$. The product estimate in the regime
$r>m/p+1$ thus gives a polynomial in the $W^{r,p}(M,f^*TN)$ norms of the sections, without
a term that would require integrating a constant over $M$. Finite multiplicity
of the cover allows the $p$ powers of the local bounds to be summed
with a common constant. This gives $\widehat P_r$ and $\widehat Q_r$.

\end{proof}

\begin{proposition}[Uniformity with fixed control data]
\label{prop:uniformidades-mapeos-cotas-fijas}
Let $(M,\mathbf{g})$ and $(N,\mathbf{h}_N)$ be complete Riemannian
manifolds without boundary satisfying $(I)$ and $(B_k)$ to the required
order. For a fixed profile $\mathbf{B}$, the families
$\mathcal V^{b,r}_\varepsilon$ restricted to
$C_{b,\mathbf{B}}^{\infty,r}(M,N)$ form a base for a separated uniformity.
The same holds in the Sobolev regime for
$\mathcal V^{p,r}_{\varepsilon,\rho}$ after choosing
$\varepsilon<\rho/C_S(\mathbf{B})$.
\end{proposition}

\begin{proof}
The diagonal belongs to every entourage, and monotonicity in
$\varepsilon$ controls intersections. Part (i) of the lemma allows us,
for each $\varepsilon>0$, to choose $\delta>0$ such that
$(\mathcal V^{b,r}_\delta)^{-1}\subseteq
\mathcal V^{b,r}_\varepsilon$; part (ii) allows $\delta$ to be decreased to
obtain
$\mathcal V^{b,r}_\delta\circ\mathcal V^{b,r}_\delta
\subseteq\mathcal V^{b,r}_\varepsilon$. The constants are uniform because
$\mathbf{B}$ is fixed. Separation follows from pointwise control. The
same argument, using the Sobolev estimates and the condition
$C_S(\mathbf{B})\delta<\rho$, proves the second assertion.
\end{proof}

\begin{theorem}[Metrization in a class with common bounds]
\label{teo:uniformidades-mapeos-eichhorn}
Let $(M,\mathbf g)$ and $(N,\mathbf h_N)$ be without boundary, and of
bounded geometry of order $k$. For $q\in\{1,\ldots,k\}$, fix a profile
$\mathbf B\in(0,+\infty)^q$. The sets
\eqref{eq:entourage-mapeos-eichhorn-b}, with $r=q$, restricted to
$C_{b,\mathbf B}^{\infty,q}(M,N)$, generate a separated uniformity with a countable
base. If $1<p<\infty$ and $k\geq r>m/p+1$, the same holds for
\eqref{eq:entourage-mapeos-eichhorn-pr} on
$C_{b,\mathbf B}^{\infty,r}(M,N)$. Both uniformities are metrizable and are
unchanged when the common normal radius is decreased.
\end{theorem}

\begin{proof}
Proposition~\ref{prop:uniformidades-mapeos-cotas-fijas} verifies the
uniformity axioms. Positive rational radii suffice to give
a countable base. Separation and metrization follow from
Proposition~\ref{prop:metrizacion-cadenas-uniformidad}. If the normal
radius is decreased, Sobolev embedding with a common constant shows that
entourages of sufficiently small radius in the original base belong to
the new base. The reverse inclusion is immediate; the bases are cofinal.
\end{proof}

\begin{example}[Individual bounds do not give a uniform radius]
\label{ej:radio-uniforme-mapeos-sin-perfil}
The need to fix the bounds can be seen directly. Take
$M=S^1\times\mathbb R$ and $N=S^2\times\mathbb R$ with their product metrics,
and, for $n\in\mathbb N$, define
\[
 f_n(\theta,t)=(\cos(n\theta),\sin(n\theta),0,t),
 \qquad \theta\in\mathbb R/(2\pi\mathbb Z),\quad t\in\mathbb R.
\]
Let $\chi\in C_c^\infty(\mathbb R)$ equal one on $[-1,1]$, with
$0\leq\chi\leq1$, and fix $0<\alpha<\rho$. The field
$\mathbf X_n=\alpha\chi(t)(0,0,1,0)$ along $f_n$ has $C_b^r(M,f_n^*TN)$ and $W^{r,p}(M,f_n^*TN)$ norms
equal to a constant times $\alpha$, independently
of $n$: the vertical vector of the sphere is parallel along the equator,
and only the factors $\chi$ are differentiated. Writing $a(t)=\alpha\chi(t)$,
\[
 h_n=\operatorname{Exp}_{f_n}\mathbf X_n
   =(\cos a\cos(n\theta),\cos a\sin(n\theta),\sin a,t).
\]
The vector giving the inverse comparison is
\[
 \mathbf Y_n=\log_{h_n}f_n
 =a(\sin a\cos(n\theta),\sin a\sin(n\theta),-\cos a,0).
\]
On $S^1\times[-1,1]$, we have
$|\nabla_{\partial_\theta}\mathbf Y_n|=n\alpha\sin\alpha$. Consequently,
the direct comparison can be arbitrarily small with a radius
independent of $n$, whereas the inverse is not bounded with that radius.
Each $f_n$ has bounded derivatives, but the family has no common profile.
Constant-radius exponential entourages on the union of all
classes do not by themselves satisfy the inversion axiom. The local
construction of mapping manifolds used below retains the
bounds of the center of each chart; see also the chart construction in
\cite[\S5]{EichhornMaps1993}.
\end{example}

\begin{remark}[Global uniformity and a local atlas]
\label{obs:cotejo-eichhorn-radio-global}
In the definition preceding Theorem~5.2 in Chapter III of
\cite{Eichhorn2007}, constant-radius entourages are considered on all
maps in the bounded-derivative class. The preceding example explains
why we do not use that assertion without fixing common bounds.
The same issue affects the Sobolev entourages in Proposition~6.1 of that
chapter: the embedding constant may depend on the center, as is made
explicit there by writing $D=D(f)$.

This does not invalidate local exponential charts. In an atlas, the
radius and constants of a chart may depend on its center; in a base
of entourages of a uniformity, the radius verifying an axiom must work
for every pair in the class under consideration. These are different
quantifiers. Moreover, a class with common bounds is not generally
open in a mapping space. Its uniform completion is not automatically
a Banach manifold. This is why the metrization in
Theorem~\ref{teo:uniformidades-mapeos-eichhorn} is kept separate from
the atlas construction in
Theorem~\ref{teo:eichhorn-variedades-mapeos-abiertas}.
\end{remark}

In the following construction, write $r$ for the order of either
norm. Fix a class with profile $\mathbf B$ and one of the two uniformities
in the theorem. Abbreviate its exponential entourages by
$\mathcal V_\varepsilon$. The inversion and composition estimates
allow us to choose recursively a decreasing rational sequence
$\varepsilon_n\to 0^{+}$ such that the symmetric entourages
\[
\mathcal W_n
:=\mathcal V_{\varepsilon_n}\cap
\mathcal V_{\varepsilon_n}^{-1}
\]
are cofinal and satisfy
$\mathcal W_{n+1}\circ\mathcal W_{n+1}\circ\mathcal W_{n+1}
\subseteq\mathcal W_n$. Define
\[
\omega(f,h):=
\inf\bigl(\{2^{-n}\mid n\in\mathbb N_0,\ (f,h)\in\mathcal W_n\}
\cup\{+\infty\}\bigr)
\]
and, for a chain $f=f_0,f_1,\ldots,f_N=h$ whose vertices
$f_i$, with $i\in\{0,\ldots,N\}$, belong to the fixed class, set
\begin{equation}
\label{eq:metrica-cadenas-eichhorn}
d_{\mathrm{cad}}(f,h)
:=
\inf\left\{
\sum_{i=0}^{N-1}\omega(f_i,f_{i+1})
\middle|
N\in\mathbb N,\ f=f_0,\ f_N=h
\right\},
\end{equation},
with value $+\infty$ when there is no chain of finite cost.
Proposition~\ref{prop:metrizacion-cadenas-uniformidad} gives
\[
\mathcal W_n
\subseteq
\{(f,h)\mid d_{\mathrm{cad}}(f,h)\leq2^{-n}\}
\quad\text{and}\quad
\{(f,h)\mid d_{\mathrm{cad}}(f,h)<2^{-n}\}
\subseteq\mathcal W_{n-1},
\qquad n\geq1.
\]
Thus \eqref{eq:metrica-cadenas-eichhorn} generates the uniformity on each
finite-distance component, and its truncation at one metrizes it on the entire
space. The construction does not presuppose that an entire component fits inside
a single exponential chart.

\begin{definition}[Slice, equivariant tube, and orbit types]
\label{def:rebanada-tubo-cuasiestratificacion-cap24}
\index{slice}
\index{equivariant tube}
Let $G$ be a Banach Lie group acting smoothly on a Banach manifold
$\mathcal X$, and let $x\in\mathcal X$. Suppose the
stabilizer $G_x$ is a closed split Banach Lie subgroup,
so that the associated quotient below has its manifold structure. A \textit{slice} at $x$ is
a submanifold $S\subseteq\mathcal X$ containing $x$, invariant under
the stabilizer $G_x$, and such that the map
\[
 G\times_{G_x}S\longrightarrow\mathcal X,
 \qquad [g,s]\longmapsto g\cdot s,
\]
is a diffeomorphism onto a $G$--invariant neighborhood of the orbit $G\cdot x$.
Here $G\times_{G_x}S$ is the quotient of $G\times S$ by
$(g,s)\sim(gh^{-1},h\cdot s)$, with $h\in G_x$. The image of this
diffeomorphism is called an \textit{equivariant tube} of the orbit.

We will also use the topological version of this notion: if $G$ is a
topological group acting continuously, the same associated map
is required to be a homeomorphism onto an invariant open set, giving a
\textit{topological slice}. Differentiable regularity of $S$ or the
tube is added as a separate conclusion. This distinction allows slices
to be formulated for Sobolev diffeomorphisms, whose operations lose
derivatives and which do not form a Banach Lie group at a fixed level.

A \textit{quasistratification by orbit types} is a partition
$\mathcal X=\coprod_{(H)}\mathcal X_{(H)}$, where
$\mathcal X_{(H)}:=\{y\in\mathcal X\mid G_y\text{ is conjugate to }H\}$,
satisfying the following conditions: each $\mathcal X_{(H)}$ is a
locally closed submanifold; slices provide local models
for the partition; and the frontier condition holds:
\[
 \mathcal X_{(K)}\cap\overline{\mathcal X_{(H)}}\neq\varnothing
 \quad\Longrightarrow\quad
 H\text{ is conjugate to a subgroup of }K.
\]
The family of types need not be countable or locally finite.
Regularity, second countability of the strata, and second countability
of the quotient are additional properties.
\end{definition}

The product in $H^1((0,1))$ was established with an explicit constant
in Lemma~\ref{lem:algebra-h1-intervalo}. On a base of dimension $m$,
products of sections in the $H^r$ scale require simultaneous control
of the order of the multiplier and that of the section on which
it acts.

\begin{lemma}[Products, derivatives, and approximation in the Hilbert scale]
\label{lem:calculo-sobolev-compacto-aplicaciones}
Let $M$ be a smooth compact manifold of dimension $m$, with or without boundary,
and let $\mathbf E,\mathbf F,\mathbf G\to M$ be smooth vector bundles of
finite rank, equipped with smooth bundle metrics and connections. The
following properties hold.
\begin{enumerate}[label=(\alph*)]
\item For $a>m/2$, $0\leq r\leq a$, and a smooth homomorphism
$B:\mathbf E\otimes\mathbf F\to\mathbf G$, there exists $C>0$ such that
\begin{equation}
\label{eq:producto-hilbert-haces-aplicaciones}
 \|B(\mathbf v\otimes\mathbf z)\|_{H^r(M,\mathbf G)}
 \leq C\|\mathbf v\|_{H^a(M,\mathbf E)}
          \|\mathbf z\|_{H^r(M,\mathbf F)}.
\end{equation}
\item For $r\geq0$, the covariant derivative is a continuous linear operator
\[
 \nabla:H^r(M,\mathbf E)\longrightarrow
 H^{r-1}(M,T^*M\otimes\mathbf E).
\]
In particular, if $s>m/2+1$, contraction of this derivative satisfies
\begin{equation}
\label{eq:estimacion-transporte-sobolev-aplicaciones}
 \|\nabla_{\mathbf v}\mathbf z\|_{H^{s-1}(M,TM)}
 \leq C\|\mathbf v\|_{H^s(M,TM)}\|\mathbf z\|_{H^s(M,TM)}.
\end{equation}
\item For every $r\geq0$, $\Gamma(\mathbf E)$ is dense in
$H^r(M,\mathbf E)$. If $r_0,r_1\geq0$ and $0<\theta<1$, then
\[
 [H^{r_0}(M,\mathbf E),H^{r_1}(M,\mathbf E)]_\theta
 =H^{(1-\theta)r_0+\theta r_1}(M,\mathbf E)
\]
with equivalent norms.
\end{enumerate}
\end{lemma}

\begin{proof}
Fix a finite cover by interior and boundary charts, smooth
frames, and a subordinate partition of unity $(\rho_i)_{i=1}^N$.
Choose cutoff functions $\chi_i$ supported in the same charts, with
$\chi_i=1$ on a neighborhood of $\operatorname{supp}\rho_i$. Local norm
equivalence for bundles and the constructions in
Theorems~\ref{teo:extension-universal-lipschitz-fraccionaria} and
\ref{teo:traza-haz-frontera-geometria-acotada} allow
$\rho_i\mathbf v$ to be represented continuously by functions in
$H^r(\mathbb R^m,\mathbb K^{\operatorname{rank}\mathbf E})$:
in boundary charts use the half-space extension, and in
interior charts extend by zero after multiplication by the cutoff function.
To return to $M$, restrict to the chart and multiply by $\chi_i$.
The sum of these operations recovers $\mathbf v$.

First prove the scalar product estimate on $\mathbb R^m$. For $a>m/2$,
Cauchy--Schwarz gives
\[
 \|\widehat v\|_{L^1(\mathbb R^m)}
 \leq\left(\int_{\mathbb R^m}\langle\xi\rangle^{-2a}\,d\xi\right)^{1/2}
      \|v\|_{H^a(\mathbb R^m)}.
\]
The integral is finite. The inequality
$\langle\xi\rangle^a\leq C_a(
\langle\xi-\zeta\rangle^a+\langle\zeta\rangle^a)$, the Fourier formula
for products, and Young's convolution inequality imply, initially
for Schwartz functions,
\[
 \|vz\|_{H^a(\mathbb R^m)}
 \leq C_a\bigl(
 \|v\|_{H^a(\mathbb R^m)}\|\widehat z\|_{L^1(\mathbb R^m)}
 +\|\widehat v\|_{L^1(\mathbb R^m)}\|z\|_{H^a(\mathbb R^m)}\bigr)
 \leq C\|v\|_{H^a(\mathbb R^m)}\|z\|_{H^a(\mathbb R^m)}.
\]
Density of Schwartz functions in
Proposition~\ref{prop:regularidad-intermedia-propiedades-escala-hilbertiana}
extends this product: if $v_j\to v$ and $z_j\to z$ in
$H^a(\mathbb R^m)$, applying the inequality to
$(v_j-v_k)z_j+v_k(z_j-z_k)$ proves that $(v_jz_j)$ is Cauchy.
Embedding into $C^0(\mathbb R^m)$ identifies its limit with $vz$.
On the other hand,
$\|vz\|_{L^2(\mathbb R^m)}\leq
\|v\|_{L^\infty(\mathbb R^m)}\|z\|_{L^2(\mathbb R^m)}$.
Interpolate these two bounds for the linear operator $z\mapsto vz$,
holding $v$ fixed;
Theorem~\ref{teo:interpolacion-escalas-BF}, with $p=q=2$, gives the bound
in $H^r(\mathbb R^m)$ for $0\leq r\leq a$.
Apply it componentwise and to the localized smooth coefficients
of $B$. The finitely many localization operations above prove (a).

In a local frame,
$\nabla_j\mathbf v$ has components
$\partial_jv^b+\displaystyle\sum_{c=1}^{\operatorname{rank}\mathbf E}\Gamma^b_{jc}v^c$.
Euclidean differentiation from the same proposition and continuity of smooth
multipliers prove (b). For the final estimate, take
$a=r=s-1>m/2$ in (a), use
$H^s(M,TM)\hookrightarrow H^{s-1}(M,TM)$, and contract
$\mathbf v\otimes\nabla\mathbf z$.

For (c), represent each $\rho_i\mathbf v$ by an extension
$V_i$ as at the beginning of the proof. Take
$V_{i,j}\in\mathcal S(\mathbb R^m,\mathbb K^{\operatorname{rank}\mathbf E})$
with $V_{i,j}\to V_i$ in the indicated space $H^r(\mathbb R^m,\mathbb K^{\operatorname{rank}\mathbf E})$. Restricting,
reconstructing in the frames, and multiplying by $\chi_i$ gives smooth
sections $\mathbf v_j$ on all of $M$, including the boundary charts, and
\[
 \|\mathbf v_j-\mathbf v\|_{H^r(M,\mathbf E)}
 \leq C\sum_{i=1}^N
 \|V_{i,j}-V_i\|_{H^r(\mathbb R^m,\mathbb K^{\operatorname{rank}\mathbf E})}
 \longrightarrow0.
\]
Finally, the localization and reconstruction operators can be
chosen simultaneously for $r_0$ and $r_1$: the Euclidean extension
operator is universal. They are therefore a coretraction and a
retraction between the section spaces and finite products of
Euclidean spaces. Interpolating both identities through
Theorem~\ref{teo:interpolacion-escalas-BF} proves the equality in (c).
\end{proof}

\section{Sobolev mapping manifolds}
\label{sec:variedades-mapeos-sobolev}

Let $(M,\mathbf{g})$ and $(N,\mathbf{h}_N)$, of dimensions $m$ and $n$, respectively, be smooth finite-dimensional
Riemannian manifolds, with $N$ without boundary; $M$ may
have boundary. Let
$1<p<\infty$, $r\in\mathbb N$. The basic condition
\begin{equation}
\label{eq:umbral-continuidad-mapeos-sobolev}
rp>m
\end{equation}
ensures that $W^{r,p}(M,N)$ classes have continuous representatives. This
continuity is indispensable: without it, one cannot require two maps to
take values in a common normal neighborhood pointwise.

\begin{definition}[Sobolev maps]
\label{def:mapeos-sobolev-variedades}
\index{Sobolev mapping space}
Let $M$ be a smooth compact manifold with or without boundary and let $N$ be a smooth
manifold without boundary. A continuous map $f\colon M\to N$ has class $W^{r,p}(M,N)$ if it admits
finite covers by charts $(U,\phi)$ on the base and $(V,\psi)$ on the target,
with $f(U)\subseteq V$, such that
$\psi\circ f\circ\phi^{-1}\in W^{r,p}(\phi(U),\mathbb R^n)$.
For a boundary chart, we use the interior of the half-ball and regularity up to its flat face.
The resulting space is denoted by $W^{r,p}(M,N)$.
\end{definition}

If $f$ is only Sobolev, $f^*TN$ is initially a continuous bundle. The notation
$W^{r,p}(M,f^*TN)$ is defined using trivializations induced by charts
of $N$: their transition matrices are smooth superpositions of $f$ and
belong locally to $W^{r,p}(\phi(U),\mathbb R^{n\times n})$. Since $rp>m$, multiplication by these
matrices and their inverses is continuous in $W^{r,p}(\phi(U),\mathbb R^{n\times n})$. A finite partition
of unity thus gives equivalent norms on the section space.
No smoothness is attributed to $f^*TN$ before its construction. If the center $f$ is
smooth, this definition agrees with the covariant norm in the chapter on
Sobolev spaces on bundles, by Lemma~\ref{lem:local-expression-higher-order}.

Fix $f\in W^{r,p}(M,N)$ with $M$ compact. The map
\[
 \mathcal E\colon TN\supset\mathcal O\longrightarrow N\times N,
 \qquad v\longmapsto\bigl(\pi(v),\exp_{\pi(v)}v\bigr),
\]
is a local diffeomorphism along the zero section. Since $f(M)$ is
compact, there exist a neighborhood $\mathcal O_f$ of the zero section over $f(M)$
and a number $\rho_f>0$ such that $|v|_{\mathbf{h}_N}<\rho_f$ implies
$v\in\mathcal O_f$ and $\mathcal E$ is injective on this set. By
Sobolev embedding,
\[
\mathcal V_f:=\{\mathbf{X}\in W^{r,p}(M,f^*TN)\mid
\|\mathbf{X}\|_{C^0(M,f^*TN;f^*\mathbf{h}_N)}<\rho_f\}
\]
is open. Define
\begin{equation}
\label{eq:carta-exponencial-mapeos-sobolev}
\operatorname{Exp}_f(\mathbf{X})(x):=\exp^N_{f(x)}\mathbf{X}(x).
\end{equation}

The exponential, logarithm, reciprocal, and inverse square root in
Example~\ref{ej:funciones-elementales-algebras} are smooth superposition
operators on the algebras considered there. Constructing geometric
charts requires this property for arbitrary smooth maps between
finite-dimensional spaces.

\begin{lemma}[Superposition in the algebra regime]
\label{lem:superposicion-sobolev-geometrica}
Let $M$ be a smooth compact manifold with or without boundary, let
$U\subseteq\mathbb R^a$ be open, and let
$F\colon U\longrightarrow\mathbb R^b$ be smooth. If $rp>m$, the superposition map
$u\mapsto F\circ u$ is smooth on the open set
$\{u\in W^{r,p}(M,\mathbb R^a)\mid u(M)\subseteq U\}$. On each neighborhood whose
values lie in a compact set $K\Subset U$, the estimates depend only on
the geometry of $M$, $r$, $p$, a bound on the $W^{r,p}(M,\mathbb R^a)$ norm of the
arguments, and bounds on the derivatives of $F$ on a neighborhood of $K$.
Its differential is
\[
D(F_*)(u)v=(DF\circ u)v.
\]
\end{lemma}

\begin{proof}
Since $rp>m$, Sobolev embedding gives each class a continuous
representative; below we identify the class with this representative. Set
\[
 \mathcal U:=\{u\in W^{r,p}(M,\mathbb R^a)\mid u(M)\subseteq U\}.
\]
First check that $\mathcal U$ is open. If $u\in\mathcal U$,
compactness of $M$ and continuity of $u$ imply that $u(M)$ is a compact set
contained in $U$. If $U\neq\mathbb R^a$, set
$d_u:=\operatorname{dist}(u(M),\mathbb R^a\setminus U)>0$; if
$U=\mathbb R^a$, simply fix $d_u:=1$. Let
$C_{\mathrm{emb}}>0$ be a constant for the embedding
$W^{r,p}(M,\mathbb R^a)\hookrightarrow C^0(M,\mathbb R^a)$. If
$\|v\|_{W^{r,p}(M,\mathbb R^a)}
<\frac{d_u}{2C_{\mathrm{emb}}}$, then
$\|v\|_{C^0(M,\mathbb R^a)}<\frac{d_u}{2}$ and
$(u+v)(M)\subseteq U$. Thus $\mathcal U$ is open.

Now fix a bounded subset
$\mathcal B\subseteq\mathcal U$ whose images lie in a common
compact set $K\Subset U$. Choose another compact set $K_1\Subset U$ whose
interior contains $K$. Decreasing $\delta>0$, the preceding embedding ensures
that, if $u\in\mathcal B$ and
$\|v\|_{W^{r,p}(M,\mathbb R^a)}<\delta$, then
$(u+tv)(M)\subseteq K_1$ for every $t\in[0,1]$.

For a smooth function $G\colon U\longrightarrow\mathbb R^c$, the
Faà di Bruno formula in Theorem~\ref{faa di bruno multivariable}, applied in charts,
expresses each derivative of order at most $r$
of $G\circ u$ as a finite sum of terms of the form
\[
 (D^\ell G\circ u)
 \bigl[D^{\alpha_1}u,\ldots,D^{\alpha_\ell}u\bigr],
 \qquad
 1\leq\ell\leq r,
 \quad
 |\alpha_1|+\cdots+|\alpha_\ell|\leq r,
 \quad
 \alpha_i\in\mathbb N_0^m,\quad |\alpha_i|\geq1\quad(1\leq i\leq\ell).
\]
To estimate each product, work on a bounded coordinate domain $\Omega\subseteq\mathbb R^m$
with fixed cutoff functions, and write $u$ for its local representation.
If $d_i=|\alpha_i|$, choose $q_i=rp/d_i$. The Sobolev embedding of
order $r-d_i$ gives
\[
 \|D^{\alpha_i}u\|_{L^{q_i}(\Omega,\mathbb R^a)}
 \leq C\|u\|_{W^{r,p}(\Omega,\mathbb R^a)}.
\]
Indeed, $1/p-(r-d_i)/m\leq d_i/(rp)$ because $rp>m$;
if $d_i=r$, simply use $q_i=p$. Since
$\displaystyle\sum_{i=1}^{\ell}1/q_i=\displaystyle\sum_{i=1}^{\ell} d_i/(rp)\leq1/p$, Hölder's inequality and the finite measure of
the region control the product in $L^p(\mathbb R^m)$.
The factor $D^\ell G\circ u$ is controlled by its supremum on $K_1$.
This argument, also applied to Leibniz's rule, proves the algebra bound
\[
 \|vw\|_{W^{r,p}(M)}
 \leq C\|v\|_{W^{r,p}(M)}\|w\|_{W^{r,p}(M)}.
\]
Summing over the multi-indices and finitely many charts gives a constant
$C=C(M,r,p,\mathcal B,K_1,G)>0$ such that
\begin{equation}
\label{eq:cota-local-superposicion-sobolev}
 \|G\circ u\|_{W^{r,p}(M,\mathbb R^c)}
 \leq C
 \qquad (u\in\mathcal B),
\end{equation}
and, if $u_0,u_1\in\mathcal B$ and the segment joining their values remains in
$K_1$,
\begin{equation}
\label{eq:lipschitz-local-superposicion-sobolev}
 \|G\circ u_1-G\circ u_0\|_{W^{r,p}(M,\mathbb R^c)}
 \leq C\|u_1-u_0\|_{W^{r,p}(M,\mathbb R^a)}.
\end{equation}
Indeed, for smooth functions the first bound follows term by term
from the preceding formula. For the second, write
\[
 G\circ u_1-G\circ u_0
 =\int_0^1
 (DG\circ(u_0+t(u_1-u_0)))(u_1-u_0)\,dt
\]
and apply the same product estimates. To extend the bounds,
let $u\in\mathcal B$ and choose, by
Meyers--Serrin~\ref{teo:meyers-serrin-haz-frontera} if there is boundary and
\ref{meyers-serrin-haz} if there is none, a sequence
$u_j\in C^\infty(M,\mathbb R^a)$ with
$u_j\to u$ in $W^{r,p}(M,\mathbb R^a)$.
Embedding into $C^0(M,\mathbb R^a)$ gives uniform convergence. After
discarding finitely many terms, the values of $u_j$ lie
in $K_1$ and their Sobolev norms have a common bound.
The segments joining $u_j(x)$ and $u_k(x)$ remain in a compact neighborhood
of $K$ contained in $U$ when the approximation is sufficiently close.
By the inequality already proved for smooth arguments,
\[
 \|G\circ u_j-G\circ u_k\|_{W^{r,p}(M,\mathbb R^c)}
 \leq C\|u_j-u_k\|_{W^{r,p}(M,\mathbb R^a)}\longrightarrow0.
\]
The limit exists by completeness and agrees with $G\circ u$ by
uniform convergence. Taking norms and limits gives
\eqref{eq:cota-local-superposicion-sobolev}. For two arguments
$u_0,u_1$, choose both approximations in the same compact set and
apply the smooth bound to their differences; their limits give
\eqref{eq:lipschitz-local-superposicion-sobolev}.
The constants depend on the stated Sobolev bounds and the
derivatives of $G$ on the chosen compact neighborhood, independently of $j$ and
$k$.

For $q\in\mathbb N_0$, define
\[
 A_q(u)[v_1,\ldots,v_q]
 :=(D^qF\circ u)[v_1,\ldots,v_q],
\]
with the convention $A_0(u)=F\circ u$. The product estimate and
\eqref{eq:cota-local-superposicion-sobolev}, applied to $G=D^qF$, show
that $A_q(u)$ is a continuous $q$--linear operator from
$W^{r,p}(M,\mathbb R^a)^q$ to $W^{r,p}(M,\mathbb R^b)$ and that
$u\mapsto A_q(u)$ is continuous in the operator norm. Moreover, for any $w_1,\ldots,w_q\in W^{r,p}(M,\mathbb R^a)$,
the fundamental theorem of calculus gives, pointwise,
\[
 \bigl(A_q(u+v)-A_q(u)-A_{q+1}(u)[v]\bigr)
 [w_1,\ldots,w_q]
 =\int_0^1(1-t)
 (D^{q+2}F\circ(u+tv))[v,v,w_1,\ldots,w_q]\,dt.
\]
Taking the operator norm and applying the preceding bounds again,
there exist $\delta,C_q>0$ such that
\[
 \|A_q(u+v)-A_q(u)-A_{q+1}(u)[v]\|_{
 \mathcal L^q(W^{r,p}(M,\mathbb R^a);W^{r,p}(M,\mathbb R^b))}
 \leq C_q\|v\|_{W^{r,p}(M,\mathbb R^a)}^2
\]
for $u\in\mathcal B$ and
$\|v\|_{W^{r,p}(M,\mathbb R^a)}<\delta$. For $q=0$, division of the remainder
by $\|v\|_{W^{r,p}(M,\mathbb R^a)}$ proves that $DF_*=A_1$. If
$D^qF_*=A_q$ has been identified, the same operator norm estimate proves that
$A_q$ is differentiable with differential $A_{q+1}$; the continuity already established
gives $F_*\in C^{q+1}$. This step, for each $q\in\mathbb N_0$, identifies
$A_q=D^qF_*$ at all orders. In particular,
\[
 D(F_*)(u)v=(DF\circ u)v,
\]
and $F_*$ is smooth on $\mathcal U$.
\end{proof}

\begin{corollary}[Superposition at real Hilbert orders]
\label{cor:superposicion-hilbert-orden-real}
The preceding lemma holds with $W^{r,p}(M,\mathbb R^a)$ replaced by $H^s(M,\mathbb R^a)$ for every real
index $s>m/2$. Consequently, the compact constructions in this chapter
using $H^s(M,\mathbf F)$, for the corresponding bundle $\mathbf F\to M$, also admit noninteger indices.
\end{corollary}

\begin{proof}
We explain the estimate that handles these indices. Localize in
a chart using cutoff functions and use the extension operator from the trace chapter
when the chart meets the boundary. It suffices to estimate superposition
on functions in $H^s(\mathbb R^m)$, multiplying the result by a
fixed cutoff function $\chi$. Around a given argument, all values lie
in a compact subset of $U$. Extend $F$ from a neighborhood of this compact set to
a smooth map on all of $\mathbb R^a$ with bounded derivatives; the extension does not
change the superposition under study.

Take the low-frequency operators $S_j$ and the blocks
$\Delta_j=S_{j+1}-S_j$, for $j\in\mathbb N_0$, of an inhomogeneous dyadic
decomposition. Convergence $S_j u\to u$ in $C^0$ and $H^s(\mathbb R^m)$ gives
the telescoping identity
\[
 F(u)=F(S_0u)+\sum_{j=0}^{\infty} a_j\Delta_j u,
 \qquad
 a_j(x)=\int_0^1 DF\bigl(S_j u(x)+t\Delta_j u(x)\bigr)\,dt.
\]
For each $N\in\mathbb N_0$, Bernstein's inequality and the chain rule give
\[
 \|D^\alpha a_j\|_\infty\leq C_N2^{j|\alpha|},
 \qquad j\in\mathbb N_0,\quad \alpha\in\mathbb N_0^m,\quad |\alpha|\leq N.
\]
Here $C_N$ depends on a bound for $\|u\|_\infty$ and the derivatives of
$F$ through order $N+1$. Applying Leibniz's rule to the product and Bernstein's inequality to the
block $\Delta_j u$ gives
$\|D^\alpha(\chi a_j\Delta_j u)\|_{L^{2}(\mathbb R^m,\mathbb R^b)}
\leq C_N2^{j|\alpha|}\|\Delta_j u\|_{L^{2}(\mathbb R^m,\mathbb R^a)}$.
If $k>j$, the vanishing moments of the kernel of $\Delta_k$ allow subtraction of the
Taylor polynomial of order $N-1$ of the product. The integral Taylor remainder,
together with the preceding bound, gives
\[
 \|\Delta_k(\chi a_j\Delta_j u)\|_{L^{2}(\mathbb R^m,\mathbb R^b)}
 \leq C_N\min\{1,2^{-N(k-j)}\}\|\Delta_j u\|_{L^{2}(\mathbb R^m,\mathbb R^a)},
 \qquad j,k\in\mathbb N_0.
\]
For $k\leq j$, we have simply used the $L^1(\mathbb R^m)$ bound on the kernel of
$\Delta_k$ and the uniform bound on $a_j$. Choose an integer $N>s$. After
multiplication by $2^{ks}$, the two regions produce the factors
$2^{-(N-s)(k-j)}$ and $2^{-s(j-k)}$, respectively. Both sequences are summable.
The convolution inequality in $\ell^2(\mathbb N_0)$ and the characterization
$H^s(\mathbb R^m)=F^s_{2,2}$ from the chapter on intermediate regularity imply
\[
 \|\chi F(u)\|_{H^s(\mathbb R^m,\mathbb R^b)}\leq C\bigl(1+\|u\|_{H^s(\mathbb R^m,\mathbb R^a)}\bigr).
\]
The low-frequency term is estimated by differentiating $\chi F(S_0u)$ through
order $N$; compact support of $\chi$ makes its $L^2(\mathbb R^m)$ norm finite.

Apply this bound also to $DF,D^2F,\ldots$. Since $H^s(\mathbb R^m)$ is an algebra,
the integral identity for $F(u+v)-F(u)$ gives local Lipschitz
continuity, and the second-order integral remainder gives differentiability
with differential $(DF\circ u)v$. For each $q\in\mathbb N_0$, applying the
same remainder to $D^qF$ identifies the differential of $A_q$ with $A_{q+1}$,
exactly as in the preceding proof. The estimates are local in the
$H^s(M,\mathbb R^b)$ norm and the cover of $M$ is finite; combining them gives the
asserted smoothness.
\end{proof}

\begin{theorem}[Exponential charts]
\label{teo:cartas-exponenciales-mapeos-sobolev}\glsadd{mapeos-sobolev}
Let $M$ be a smooth compact manifold with or without boundary and $N$ a smooth
finite-dimensional Riemannian manifold without boundary, and suppose
$rp>m$. For sufficiently small $\rho_f$,
the map \eqref{eq:carta-exponencial-mapeos-sobolev} is a diffeomorphism from
$\mathcal V_f$ onto an open subset of $W^{r,p}(M,N)$. These charts define a
Banach manifold independent of the auxiliary metric on $N$. The operator
\[
J_f:=D(\operatorname{Exp}_f^{-1})_f
\colon T_fW^{r,p}(M,N)
\longrightarrow W^{r,p}(M,f^*TN)
\]
is a continuous linear isomorphism with inverse
$D\operatorname{Exp}_f(0)$. In particular,
\begin{equation}
\label{eq:tangente-mapeos-sobolev}
T_fW^{r,p}(M,N)\cong W^{r,p}(M,f^*TN)
\end{equation}
canonically.
\end{theorem}

The construction goes back to \cite[Sections 6--8]{Eells1966} and
\cite[Chapters II--III]{Palais1968}.

\begin{proof}
Let $\mathcal W_f:=\mathcal E(\mathcal O_f)$. The inverse of
$\operatorname{Exp}_f$ on the set of $h$ for which
$(f(x),h(x))\in\mathcal W_f$ for every $x$ is
\[
 \log_f(h)(x):=\mathcal E^{-1}(f(x),h(x))
              =\exp_{f(x)}^{-1}h(x).
\]
The condition on $(f,h)$ is open in the $C^0$ topology because
$M$ is compact and the graph of $(f,h)$ is compact inside
$\mathcal W_f$; by $W^{r,p}(M,f^*TN)\hookrightarrow C^0(M,f^*TN)$, it is also open in the
Sobolev topology.

To prove smoothness, choose a finite cover of $M$ over which
$f^*TN$ is trivialized, and charts of $N$ containing the relevant compact sets.
In these coordinates, $\mathbf{X}\mapsto\operatorname{Exp}_f\mathbf{X}$ and
$h\mapsto\log_fh$ are superposition operators given by smooth maps whose
derivatives are bounded on a compact set.
Lemma~\ref{lem:superposicion-sobolev-geometrica} proves that they are smooth; a
finite partition of unity combines the local expressions.

If the charts are centered at $f$ and $g$, the coordinate change has
the intrinsic expression
\begin{equation}
\label{eq:cambio-cartas-mapeos-sobolev}
 \Theta_{g f}(\mathbf{X})(x)
 :=\exp_{g(x)}^{-1}\bigl(\exp_{f(x)}\mathbf{X}(x)\bigr).
\end{equation}
The preceding argument, applied to the smooth map
$(y,z,v)\mapsto\exp_z^{-1}(\exp_yv)$ on its domain of definition,
proves that \eqref{eq:cambio-cartas-mapeos-sobolev} is smooth. This directly verifies
compatibility of the atlas. If the auxiliary metric is changed,
the same formula with the two exponential maps shows that the identity between the
two atlases is smooth in both directions. Finally,
$D\operatorname{Exp}_f(0)\mathbf{X}=\mathbf{X}$ pointwise. Since
$\operatorname{Exp}_f$ and its inverse are smooth, their differentials are
continuous; differentiating the two inverse identities gives
$J_fD\operatorname{Exp}_f(0)=\operatorname{id}$ and
$D\operatorname{Exp}_f(0)J_f=\operatorname{id}$. This identification is
independent of the chart and gives \eqref{eq:tangente-mapeos-sobolev}.

For $M=I$ and $N=S^{d-1}$ with its usual metric, exponential charts
are compatible with the normalization charts in
Example~\ref{ej:curvas-esfera-extremos-banach}. Identify fields along
$\gamma$ with $H^1(I,\mathbb R^d)$ functions orthogonal
to $\gamma$ pointwise. On the overlap domains, the coordinate
changes between the two inverse charts are
\[
\begin{aligned}
 \Psi_\gamma(h)(t)
 &=\exp_{\gamma(t)}^{-1}
   \left(\frac{\gamma(t)+h(t)}
   {\bigl(1+|h(t)|^2\bigr)^{\frac{1}{2}}}\right),\\
 \Xi_\gamma(\mathbf V)(t)
 &=\frac{\exp_{\gamma(t)}\mathbf V(t)}
   {\gamma(t)\cdot\exp_{\gamma(t)}\mathbf V(t)}-\gamma(t).
\end{aligned}
\]
The first domain requires the normalized curve to lie in the domain
of the exponential logarithm; the second also requires strict
positivity of the denominator. These conditions are open in $C^0$
by compactness of $I$, hence in $H^1$. On a finite cover by target
charts, both expressions are superpositions of smooth maps of
the variables $\gamma(t)$ and $h(t)$ or $\mathbf V(t)$.
Lemma~\ref{lem:superposicion-sobolev-geometrica}, with $r=1$ and $p=2$,
proves their smoothness. The inverse identities for the exponential
and logarithm, and those for normalization in that example, give
$\Xi_\gamma\Psi_\gamma(h)=h$ and
$\Psi_\gamma\Xi_\gamma(\mathbf V)=\mathbf V$ on the corresponding
domains. With fixed endpoints, both transformations preserve
$h(0)=h(1)=0$ and $\mathbf V(0)=\mathbf V(1)=0$.
Thus the example and the general construction induce the same
smooth structure on the manifold of sphere-valued curves
with fixed endpoints.
\end{proof}

\subsection{Composition estimates and loss of derivatives}
\label{subsec:composicion-sobolev-prueba}

Before stating evaluation and composition properties, we prove their
underlying change-of-variable estimate. Here $\operatorname{Diff}^s(M)$
denotes the set of $C^1$ diffeomorphisms of $M$ whose expressions as
maps are $H^s(M,M)$; it carries the open mapping-chart structure, with the
appropriate boundary constraints. Differentiability of inversion is a
separate issue and is not needed for the following estimate.

\begin{lemma}[Module property in the Hilbert scale]

\label{lem:modulo-hilbert-orden-real-composicion}

If $a>m/2$ and $0\leq b\leq a$, multiplication is continuous from
$H^a(M)\times H^b(M)$ to $H^b(M)$ for compact $M$, with or without
boundary: $\|uv\|_{H^b(M)}\leq C\|u\|_{H^a(M)}\|v\|_{H^b(M)}$.

\end{lemma}

\begin{proof}

On $\mathbb R^m$, Cauchy--Schwarz and $a>m/2$ imply
$\|\widehat u\|_{L^1(\mathbb R^m)}\leq C_a\|u\|_{H^a(\mathbb R^m)}$.
Apply $\langle\xi+\eta\rangle^a\leq
C_a(\langle\xi\rangle^a+\langle\eta\rangle^a)$ to the convolution
representing $\widehat{uv}$, and then use Young to obtain

\[
 \|uv\|_{H^a(\mathbb R^m)}\leq C_a\bigl(
 \|u\|_{H^a(\mathbb R^m)}\|\widehat v\|_{L^1(\mathbb R^m)}
 +\|\widehat u\|_{L^1(\mathbb R^m)}\|v\|_{H^a(\mathbb R^m)}\bigr)
 \leq C_a'\|u\|_{H^a(\mathbb R^m)}\|v\|_{H^a(\mathbb R^m)}.
\]

Also, $\|uv\|_{L^2(\mathbb R^m)}\leq\|u\|_\infty\|v\|_{L^2(\mathbb R^m)}
\leq C_a\|u\|_{H^a(\mathbb R^m)}\|v\|_{L^2(\mathbb R^m)}$. For fixed $u$, interpolate the
linear operator $v\mapsto uv$ between $H^0(\mathbb R^m)$ and $H^a(\mathbb R^m)$.
The Hilbert-scale interpolation identity in
Theorem~\ref{teo:interpolacion-escalas-BF} gives the $H^b(\mathbb R^m)$ bound.
Finitely many cutoffs and charts, with local boundary extensions,
transport the estimate to $M$. The same argument applies to bundle
components and finite-rank tensor contractions.

\end{proof}

\begin{lemma}[Composition by a Sobolev diffeomorphism]

\label{lem:composicion-sobolev-autocontenida}

Let $M$ be compact, $m=\dim M$, $s>m/2+1$, and
$\eta\in\operatorname{Diff}^s(M)$. For $0\leq r\leq s$,
$u\mapsto u\circ\eta$ is bounded linear on $H^r(M)$.
Its operator norms are locally bounded as $\eta$ varies in
$\operatorname{Diff}^s(M)$ with the topology induced by $H^s(M,M)$.
The joint map $(u,\eta)\mapsto u\circ\eta$ is continuous from
$H^r(M)\times\operatorname{Diff}^s(M)$ to $H^r(M)$.

\end{lemma}

\begin{proof}

First work on relatively compact coordinate domains with cutoffs and
a bi-Lipschitz diffeomorphism between them. For $r=0$, change of
variables gives

\[
 \|u\circ\eta\|_{L^2(\Omega)}^2
 =\int_{\eta(\Omega)}|u(y)|^2|\det D\eta^{-1}(y)|\,dy
 \leq\|\det D\eta^{-1}\|_\infty\|u\|_{L^2(\eta(\Omega))}^2.
\]

For $0<r<1$, use the equivalent Slobodeckij norm already constructed.
If $|\eta(x)-\eta(z)|\leq L|x-z|$, then

\begin{align*}
 &\int_{\Omega}\int_{\Omega}
   \frac{|u(\eta(x))-u(\eta(z))|^2}{|x-z|^{m+2r}}\,dx\,dz\\
 &\quad\leq L^{m+2r}\|\det D\eta^{-1}\|_\infty^2
 \int_{\eta(\Omega)}\int_{\eta(\Omega)}
   \frac{|u(y)-u(w)|^2}{|y-w|^{m+2r}}\,dy\,dw.
\end{align*}

For $r=1$, also use $D(u\circ\eta)=(Du\circ\eta)D\eta$, first for
smooth $u$ and then by density. These three estimates require only
bi-Lipschitz and Jacobian bounds.

For $r>1$, proceed inductively in intervals of length one. The
derivative characterization and the Sobolev module property state that
the $H^r(\mathbb R^m)$ norm is controlled by the $L^2(\mathbb R^m)$ norm and the $H^{r-1}(\mathbb R^m)$ norms
of first derivatives, and multiplication by an $H^{s-1}(\mathbb R^m)$ element is
bounded on $H^{r-1}(\mathbb R^m)$ when $0\leq r-1\leq s-1$ and $s-1>m/2$.
The induction hypothesis therefore gives

\[
 \|(Du\circ\eta)D\eta\|_{H^{r-1}(\mathbb R^m)}
 \leq C\|Du\circ\eta\|_{H^{r-1}(\mathbb R^m)}\|D\eta\|_{H^{s-1}(\mathbb R^m)}
 \leq C_\eta\|u\|_{H^r(\mathbb R^m)}.
\]

This gives the estimate for smooth $u$; approximation in $H^r(\mathbb R^m)$ extends
it to every $u$. The limit agrees almost everywhere with composition
by the $L^2(\mathbb R^m)$ estimate. To pass to the manifold, fix finitely many
charts and cutoffs whose images under $\eta$ lie in target charts,
extend the cutoff components, and add the estimates. At the boundary,
use the half-space extensions from the Sobolev chapters; only a fixed
finite order is involved. Bounds persist as $\eta$ varies in an
$H^s(M,M)$ neighborhood: the embedding controls $C^1$, matrix inversion
bounds inverse Jacobians, and the $H^{s-1}(U,\mathbb R^{m\times m})$ norm of $D\eta$ remains
bounded.

For joint continuity, let $u_n\to u$ in $H^r(M)$ and $\eta_n\to\eta$
in $H^s(M,M)$. Local boundedness controls $(u_n-u)\circ\eta_n$.
Approximate $u$ by a smooth map $u_\varepsilon$: the two composition
errors for $u-u_\varepsilon$ are bounded by
$C\|u-u_\varepsilon\|_{H^r(M)}$ uniformly in $n$.
For fixed $u_\varepsilon$, smooth superposition gives
$u_\varepsilon\circ\eta_n\to u_\varepsilon\circ\eta$ in $H^s(M)$,
hence in $H^r(M)$. First let $n\to\infty$, then $\varepsilon\to0$.

\end{proof}

\begin{lemma}[Differentiability with controlled loss]

\label{lem:composicion-sobolev-perdida-autocontenida}

For $j\in\mathbb N_0$ and $s>m/2+1$, composition is $C^j$ as a map

\[
 H^{s+j}(M,N)\times\operatorname{Diff}^s(M)\longrightarrow H^s(M,N).
\]

In particular, $j=0$ asserts only continuity.

\end{lemma}

\begin{proof}

The assertion is local in both maps. Choose mapping charts and
finitely many base charts as in the preceding lemma. Additional
transitions are smooth superpositions, so it suffices to treat
Euclidean components. There, the candidate first derivative is

\[
 D\operatorname{Comp}(u,\eta)[h,v]
       =h\circ\eta+(Du\circ\eta)v.
\]

The preceding lemma and the $H^s(M)$ algebra show that this expression
is bounded in $(h,v)$ when $u\in H^{s+1}(M)$. The remainder of
$u\circ(\eta+v)$ after its linear term is

\[
 \int_0^1\bigl(Du\circ(\eta+tv)-Du\circ\eta\bigr)v\,dt.
\]

Work in a $C^1$ neighborhood where intermediate maps are locally
bi-Lipschitz in the chosen charts. Let $C_s$ be the $H^s(\mathbb R^m)$ multiplication constant in these charts.
Given $u,\eta$, and $\varepsilon>0$, joint continuity in the preceding
lemma provides $\delta>0$ such that
\[
 \sup_{0\leq t\leq1}\|Du\circ(\eta+tv)-Du\circ\eta\|_{H^s(\mathbb R^m)}
 \leq\varepsilon/C_s,\qquad \|v\|_{H^s(\mathbb R^m)}<\delta.
\]
Reduce $\delta$ to preserve the chosen $C^1$ neighborhood. The integral
remainder then has norm at most $\varepsilon\|v\|_{H^s(\mathbb R^m)}$. The increment in $h$ also produces
$h\circ(\eta+v)-h\circ\eta$, whose norm is at most
$C\|h\|_{H^{s+1}(\mathbb R^m)}\|v\|_{H^s(\mathbb R^m)}$ by the same integral formula.
This proves Fréchet differentiability, not merely directional
differentiability.

Operator-norm continuity of $h\mapsto h\circ\eta$ follows from the
last estimate with $\|h\|_{H^{s+1}(\mathbb R^m)}\leq1$. Continuity of the second
term follows from joint continuity of $Du\circ\eta$ in $H^s(\mathbb R^m)$ and
continuous multiplication. For higher derivatives, the exact component
formula, for $1\leq q\leq j$, is

\begin{align*}
 &D^q\operatorname{Comp}(u,\eta)
       [(h_1,v_1),\ldots,(h_q,v_q)]\\
 &\quad=(D^qu\circ\eta)[v_1,\ldots,v_q]
       +\displaystyle\sum_{i=1}^q
          (D^{q-1}h_i\circ\eta)
          [v_1,\ldots,\widehat v_i,\ldots,v_q].
\end{align*}

Inductively differentiate the preceding formula once: dependence on
$u$ is linear, so products of two $h_i$ increments do not occur; each
$\eta$ derivative adds a spatial derivative to the outer factor.
The $H^{s+j}(M)$ norm supplies every outer derivative up to order $j$
in $H^s(M)$. Multiplication in $H^s(M)$ bounds the multilinear expressions.
For continuity, the term with $D^qu$ uses joint continuity of
composition; in the terms with $h_i$, the additional derivative
available from $h_i\in H^{s+j}(M)\subseteq H^{s+q}(M)$ gives operator-norm
continuity by the integral formula already proved. Remainders are
controlled by the same formula at each step. This proves $C^j$.
Charts and cutoffs transport the estimates to manifolds; at a source
boundary, extend components locally before taking increments.

\end{proof}

\begin{proposition}[Evaluation and composition]
\label{prop:evaluacion-composicion-sobolev}
Let $M$ be a smooth compact manifold with or without boundary and let $N$ and $P$ be smooth
manifolds without boundary. Under the remaining assumptions of the preceding theorem:
\begin{enumerate}[label=(\roman*)]
\item for each $x\in M$, evaluation
$\operatorname{ev}_x(f)=f(x)$ is smooth;
\item joint evaluation
$\operatorname{ev}\colon W^{r,p}(M,N)\times M\longrightarrow N$ is continuous;
\item if $r>\frac{m}{p}+j$, joint evaluation is of class $C^j$;
\item if $H\colon N\longrightarrow P$ is smooth, the map $f\mapsto H\circ f$ is smooth on all of $W^{r,p}(M,N)$; the necessary derivative bounds
are local around the compact image of each map, not global bounds on $N$;
\item the map $(f,\eta)\mapsto f\circ\eta$ loses derivatives: in the Hilbert case,
for $s>\frac{m}{2}+1$,
\[
H^{s+j}(M,N)\times\operatorname{Diff}^{s}(M)\longrightarrow H^s(M,N),
\qquad(f,\eta)\longmapsto f\circ\eta,
\]
is of class $C^j$, whereas at the same level $H^s(M)\times
\operatorname{Diff}^s\longrightarrow H^s(M)$ only continuity is asserted in general.
\end{enumerate}
\end{proposition}

The last item is the derivative-loss composition statement used in \cite[Section 2]{EbinMarsden1970}; the preceding lemmas supply its proof here.

\begin{proof}
In exponential coordinates, evaluation has the form
$(\mathbf{X},x)\mapsto\exp_{f(x)}\mathbf{X}(x)$. The continuous embedding
$W^{r,p}(M,f^*TN)\hookrightarrow C^0(M,f^*TN)$ makes point evaluation continuous, and if
$r>\frac{m}{p}+j$, embedding $W^{r,p}(M,f^*TN)\hookrightarrow C^j(M,f^*TN)$ allows this formula
to be differentiated $j$ times with respect to $x$ and $\mathbf{X}$; each differential is a finite
sum of evaluations of derivatives of $\mathbf{X}$ multiplied by smooth
coefficients of the exponential map. This proves the first three items. The fourth is
Lemma~\ref{lem:superposicion-sobolev-geometrica}. For the last, the first
formal derivative with respect to $\eta$ contains
\[
D_\eta(f\circ\eta)\,\mathbf{Y}=(df\circ\eta)\mathbf{Y},
\]
and uses one derivative of $f$. Lemma~\ref{lem:composicion-sobolev-perdida-autocontenida} proves the remainder estimates and operator-norm continuity of all derivatives up to order $j$, giving the assertion on these exact domains.
\end{proof}

\subsection{Mappings over bounded-geometry bases without boundary}
\label{subsec:mapeos-base-no-compacta}

Without assuming compactness of the base, charts require common estimates
in the spatial variable, but their constants may depend on the center. It is useful to
construct the atlas first and keep this issue separate from
completion of a class with fixed profile.

\begin{lemma}[Mixed derivatives of exponential comparisons]
\label{lem:jets-mixtos-exponencial-orden-finito}
Let $(N,\mathbf h)$ be complete, without boundary, and of bounded
geometry of order $k\geq1$. There is $\rho>0$, with
$4\rho<\operatorname{inj}(N)$, for which the following statements hold.
The mixed covariant derivatives of the exponential map and radial
parallel transport through order $k+1$ are uniformly bounded when the
radial vector has length less than $3\rho$. These include derivatives
with respect to the initial point and the vector; the Levi--Civita
connection separates the two variables. The normal inverse
$(y,z)\mapsto\log_yz$ satisfies the same bounds for
$d_{\mathbf h}(y,z)<3\rho$.

Let $M$ be without boundary, and of bounded geometry of
order $k$, let $1<p<\infty$, and let
$k\geq r>\dim(M)/p+1$, with $k,r\in\mathbb N$.
For $f,g\in C_b^{\infty,r}(M,N)$ with
$d_{\mathbf h}(f(x),g(x))<2\rho$ at every $x\in M$, set
\[
 \Psi(x,v)=\log_{g(x)}(\exp_{f(x)}v),\qquad |v|<\rho,
 \qquad q=k+1-r.
\]
The mixed coefficients needed to differentiate this expression $r$ times
in $x$ and $q$ times in $v$ are bounded in uniform trivializations of
$f^*TN$ and $g^*TN$. Their constants may depend on the bounds for the
derivatives of $f$ and $g$ through order $r$.

For each coefficient $C$ involving a geometric derivative of order
$k+1$, there is a common constant $A$ such that
\begin{equation}\label{eq:modulo-jets-comparacion-diagonal}
 |C(x,v)-C(x,w)|
 \leq A|v-w|+
 A\bigl(d_{\mathbf h}(f(x),g(x))+|v|+|w|\bigr),
 \qquad |v|,|w|<\rho.
\end{equation}
Coefficients involving geometric derivatives only through order $k$
are uniformly Lipschitz in $v$. Over a compact subset of the base, all
the coefficients are uniformly continuous on smaller closed domains.
\end{lemma}

\begin{proof}
Write $R$ for the curvature of $N$ and
$K_j=\|\nabla^jR\|_\infty$, for $0\leq j\leq k$.
Consider a family of geodesics
$c(\boldsymbol s,t)=\exp_{y(\boldsymbol s)}(t v(\boldsymbol s))$,
where $0\leq t\leq1$. Put
$T=\partial_tc$, $J_a=\partial_{s_a}c$, and
$D_a=\nabla_{\partial_{s_a}}$. The initial data are
\[
 J_a(0)=\partial_{s_a}y,\qquad
 D_tJ_a(0)=\nabla_{\partial_{s_a}}v.
\]
These two entries vary the initial point and the vector independently.
The Jacobi equation and the commutation identities are
\[
 D_t^2J_a+R(J_a,T)T=0,
 \qquad [D_a,D_t]U=R(J_a,T)U,
 \qquad D_aT=D_tJ_a.
\]
Thus the recurrence
\eqref{eq:indexacion-exacta-fuentes-jacobi} also holds when the initial
point varies: its initial data change, but its six commutation terms do
not. Indeed, with $\mathcal LY=D_t^2Y+R(Y,T)T$, one has
\begin{align*}
 \mathcal L(D_aY)
 &=D_a(\mathcal LY)-2R(J_a,T)D_tY
   -(\nabla_TR)(J_a,T)Y-R(D_tJ_a,T)Y\\
 &\quad-(\nabla_{J_a}R)(Y,T)T
   -R(Y,D_tJ_a)T-R(Y,T)D_tJ_a.
\end{align*}
Leibniz's rule, this identity, and $D_aT=D_tJ_a$ enumerate every
term in the successive derivatives.

We specify the higher initial data. To compute a derivative at $(y,v)$,
we may take
\[
 y(\boldsymbol s)=\exp_y\left(\sum_{a=1}^{d}s_a u_a\right),\qquad
 v(\boldsymbol s)=P_{y\to y(\boldsymbol s)}
              \left(v+\sum_{a=1}^{d}s_a w_a\right),
 \qquad |u_a|+|w_a|\leq1,
\]
where $d$ is the number of differentiating directions and transport is
radial. Evaluate at $\boldsymbol s=0$. To control the data there, use
the Taylor identities in normal coordinates and a radial frame.
Here is the curvature order involved. If $A_i(z)$ is the connection
matrix in the radial frame and $\mathcal R_{ji}(z)$ its curvature in
coordinate directions, then
\[
 A_i(z)=\sum_{j=1}^{n}z^j\int_0^1t\mathcal R_{ji}(tz)\,dt,
 \qquad 1\leq i\leq n.
\]
Indeed, $\displaystyle \sum_{j=1}^{n}z^jA_j(z)=0$. Contracting the curvature
formula with $z^j$ gives the derivative of $tA_i(tz)$, which we
integrate from zero to one. A derivative of $A_i$ of order $b\geq1$
at zero contains curvature derivatives through order $b-1$.
Converting partial to covariant derivatives only adds products of lower
orders, so induction on $b$ proves this assertion. The Jacobi equation
for the coordinate frame, with data $(0,u)$, likewise expresses the
derivatives of the Christoffel symbols of order $b$ at zero using
curvature through order $b-1$. Both types of coefficients vanish at zero.

A covariant derivative of $y(\boldsymbol s)$ of order $\ell$ uses
Christoffel symbols through order $\ell-2$; one of
$v(\boldsymbol s)$ uses connection matrices through order $\ell-1$.
Thus their data of order $\ell\geq2$ require curvature only through
orders $\max\{\ell-3,0\}$ and $\ell-2$, respectively. The first
order data were written above. Terms introduced by commutation with
$D_t$ come from $R(J_a,T)$ and obey the same order bounds.
Consequently, the initial data through order $k+1$ are bounded by
$K_0,\ldots,K_{k-1}$ and depend affinely on $v$.

Now argue by induction on $\ell$. The unknown representing a derivative
of order $\ell$ satisfies an inhomogeneous Jacobi equation. Its source
contains fields of orders less than $\ell$ and curvature derivatives
through order $\ell-1$, by the recurrence above. In a parallel frame
the corresponding integral equation is
\[
 Y(t)=Y(0)+tD_tY(0)+\int_0^t(t-\tau)
       \bigl(Q(\tau)-R(Y(\tau),T(\tau))T(\tau)\bigr)\,d\tau.
\]
Gr\"onwall's inequality for the system $(Y,D_tY)$ gives
\[
 \sup_{t\in[0,1]}\bigl(|Y(t)|+|D_tY(t)|\bigr)
 \leq C\left(|Y(0)|+|D_tY(0)|+\int_0^1|Q(t)|\,dt\right).
\]
The constant depends on $K_0$ and the fixed radius. The initial data
and the induction hypothesis bound the right hand side. For parallel
transport use $D_tU=0$ and its differentiated identity
$D_tD_aU=-R(J_a,T)U$. Its derivative of order $\ell$ uses curvature
through order $\ell-1$. This proves the bounds through order $k+1$.

We retain additional information at the highest order. In the source of
order $\ell$, terms containing $\nabla^{\ell-1}R$ retain a factor
$T$. To check this, start with the six terms above. The two containing
$\nabla R$ have radial factors. Raising the curvature order at every
step requires differentiating the curvature factor without consuming
all the radial factors. Differentiating one of the latter instead leaves
the curvature order at least one below the maximum. The transport
equation has the same property, starting with $R(J_a,T)U$.
Since $|T(t)|=|v|$, the part of the source of order $k+1$ containing
$\nabla^kR$ has norm at most $C|v|$.

Separate this part of the source and solve the two linear equations with
the same principal operator. The solution it produces with zero initial
data has norm at most $C|v|$. The remaining solution has a uniformly
bounded first vertical derivative: differentiating its source once
introduces only $\nabla^jR$ with $j\leq k$ and field derivatives of
orders at most $k+1$, already bounded. The same holds when differentiating
the principal operator and the initial data. Gr\"onwall's inequality
again gives the bound. Hence every geometric derivative of order $k+1$
is the sum of a part uniformly Lipschitz in $v$ and a part bounded by
$C|v|$. The highest order of curvature differentiation in these bounds is $k$.

Reduce $\rho$ so that the vertical differential of the exponential map,
transported back to the initial point, differs from the identity in norm
by less than $1/2$ when $|v|<3\rho$. The integral Jacobi equation with
data $(0,w)$ gives this uniformly in $y$. Its inverse has norm at most
$2$. Differentiating $\exp_y(\log_yz)=z$ and multiplying by this
inverse isolates the highest derivative of the logarithm at each order.
The chain rule produces only products of controlled derivatives, and a
derivative of order $k+1$ occurs linearly. Thus the logarithm inherits
the same decomposition, with a remainder bounded by
$C\,d_{\mathbf h}(y,z)$.

Finally apply the chain rule to $\log_{g(x)}(\exp_{f(x)}v)$.
Each derivative in $x$ differentiates a center, an entry, or a coefficient,
and uses the centers only through order $r$.
Conversions between covariant derivatives and components in synchronous
frames use connection coefficients through order $r-1$, as in the
expansions of Proposition~\ref{prop:covariantes-localizacion-orden-entero}.
A geometric derivative of order $k+1$ occurs at most once in each
product. Its remainder is bounded by
$C(|v|+d_{\mathbf h}(g(x),\exp_{f(x)}v))$, which the triangle
inequality bounds by $C(2|v|+d_{\mathbf h}(f(x),g(x)))$.
Subtracting the expressions at $v$ and $w$ proves
\eqref{eq:modulo-jets-comparacion-diagonal}. At lower orders,
differentiate once more in $v$ and use the bounds through order $k+1$.
Over a compact subset of $M$, the last assertion follows from
smoothness and compactness of the closed bundle disks in use.
\end{proof}

Fix $k\geq r>m/p+1$ and a normal radius $\rho$ as in
Lemma~\ref{lem:jets-mixtos-exponencial-orden-finito}. For each smooth center
$f\in C_b^{\infty,r}(M,N)$, set
\[
 \mathcal V_f=\{\mathbf X\in W^{r,p}(M,f^*TN)\mid\|\mathbf X\|_\infty<\rho\},
 \qquad \mathcal U_f=\operatorname{Exp}_f(\mathcal V_f),
 \qquad \Omega^{p,r}(M,N)=\bigcup_{f\in C_b^{\infty,r}(M,N)}\mathcal U_f.
\]
The elements of this union are actual maps; if two exponential
images agree pointwise, they represent the same element. Give
the union the topology of these charts. Compatibility and openness of
overlaps are checked in the following proof.

\begin{theorem}[Sobolev charts in bounded geometry without boundary]
\label{teo:eichhorn-variedades-mapeos-abiertas}
Let $(M,\mathbf g)$ and $(N,\mathbf h)$, of dimensions $m$ and $n$, respectively, be smooth, without boundary,
and of bounded geometry of order $k$. Let $1<p<\infty$ and
$k\geq r>m/p+1$, with $k,r\in\mathbb N$. The preceding charts equip
$\Omega^{p,r}(M,N)$ with a Banach manifold structure of class
$C^{k+1-r}$, and a Hilbert manifold structure if $p=2$. In the case of bounded geometry of all
orders, the atlas is smooth. If $\Phi_f(\mathbf Y)=\operatorname{Exp}_f\mathbf Y$
is a chart centered at $f$, then
\[
 J_f:=D(\Phi_f^{-1})_f:\ T_f\Omega^{p,r}(M,N)
       \longrightarrow W^{r,p}(M,f^*TN)
\]
is linear and bicontinuous, with inverse $D\Phi_f(0)$.
\end{theorem}

\begin{proof}
The embedding $W^{r,p}(M,f^*TN)\hookrightarrow C_b^1(M,f^*TN)$ makes $\mathcal V_f$
open. The choice of normal radius makes $\operatorname{Exp}_f$
injective there. We use the bounds of
Lemma~\ref{lem:jets-mixtos-exponencial-orden-finito}, reducing the common
radius if necessary.

Consider a point in an overlap,
\[
 h=\operatorname{Exp}_f\mathbf X_0
   =\operatorname{Exp}_g\mathbf Y_0,
 \qquad \mathbf X_0\in\mathcal V_f,\quad
 \mathbf Y_0\in\mathcal V_g.
\]
Then $d_{\mathbf h}(f(x),g(x))\leq
|\mathbf X_0(x)|+|\mathbf Y_0(x)|<2\rho$ for all $x\in M$.
Where defined, the coordinate change is
\[
 \mathcal T(\mathbf X)(x)=\Psi(x,\mathbf X(x)),\qquad
 \Psi(x,v)=\log_{g(x)}(\exp_{f(x)}v).
\]
On a sufficiently small ball around $\mathbf X_0$ the pointwise
expression stays in the normal domain of the lemma. We prove that it
takes values in $W^{r,p}(M,g^*TN)$ and is of class $C^q$, where
$q=k+1-r\geq1$.

For $1\leq a\leq q$, define the candidate derivatives by
\begin{equation}\label{eq:derivadas-transicion-sobolev-explicitas}
 A_a(\mathbf X)[\mathbf H_1,\ldots,\mathbf H_a](x)
 =D_v^a\Psi(x,\mathbf X(x))
          [\mathbf H_1(x),\ldots,\mathbf H_a(x)],
 \qquad x\in M.
\end{equation}
Here $\mathbf H_i\in W^{r,p}(M,f^*TN)$ for $1\leq i\leq a$.
Work in uniform trivializations of the two bundles. For
$\nu\in\mathbb N_0^m$ with $|\nu|\leq r$, applying
$\partial^\nu$ to the right hand side gives a finite sum of contractions
of the form
\begin{equation}\label{eq:terminos-transicion-indices-sobolev}
 (\partial_x^\alpha D_v^{a+b}\Psi)(x,\mathbf X(x))
 \left[\partial^{\beta_1}\mathbf X,\ldots,
       \partial^{\beta_b}\mathbf X,
       \partial^{\gamma_1}\mathbf H_1,\ldots,
       \partial^{\gamma_a}\mathbf H_a\right],
\end{equation}
where $0\leq b\leq|\nu|$,
$\beta_j\in\mathbb N_0^m\setminus\{0\}$ for $1\leq j\leq b$,
$\alpha,\gamma_i\in\mathbb N_0^m$ for $1\leq i\leq a$, and
\[
 \alpha+\sum_{j=1}^{b}\beta_j+\sum_{i=1}^{a}\gamma_i=\nu.
\]
The integer coefficients of this sum follow from Leibniz's rule and the
chain rule; for $b=0$ there are no entries coming from $\mathbf X$.
The geometric order of each coefficient satisfies
\[
 |\alpha|+a+b\leq |\nu|+a\leq r+q=k+1.
\]
The preceding lemma and the center bounds therefore bound all these
coefficients uniformly.

We spell out the product estimate. If the derivative orders of the
factors are $d_1,\ldots,d_s$, then $0\leq d_i\leq r$ and
$\displaystyle \sum_{i=1}^{s}d_i\leq r$.
Write the factors as $\partial^{\delta_i}u_i$, with
$\delta_i\in\mathbb N_0^m$ and $|\delta_i|=d_i$, for $1\leq i\leq s$.
Choose $\theta_i\in[d_i/r,1]$ with
$\displaystyle \sum_{i=1}^{s}\theta_i=1$, and put $1/p_i=\theta_i/p$, interpreting
$p_i=\infty$ when $\theta_i=0$. The uniform Sobolev embeddings give
\[
 \|\partial^{\delta_i}u_i\|_{L^{p_i}(\mathbb R^m)}
 \leq C\|u_i\|_{W^{r,p}(\mathbb R^m)},\qquad 1\leq i\leq s.
\]
For $d_i=r$, this is the $L^p(\mathbb R^m)$ bound for that derivative. If $d_i<r$,
the embedding order inequality is strict because $r>m/p$; when
$p_i=\infty$, use $d_i=0$ and the uniform embedding.
H\"older's inequality applies since $\displaystyle \sum_{i=1}^{s}1/p_i=1/p$.
Localization and bounded multiplicity of the cover now give
\begin{equation}\label{eq:cota-operador-transicion-sobolev}
 \|A_a(\mathbf X)[\mathbf H_1,\ldots,\mathbf H_a]\|_{W^{r,p}(M,g^*TN)}
 \leq C_a(1+\|\mathbf X\|_{W^{r,p}(M,f^*TN)})^r
       \prod_{i=1}^{a}\|\mathbf H_i\|_{W^{r,p}(M,f^*TN)}.
\end{equation}
There is always a section factor, since $a\geq1$. No term requires a
constant function to belong to $L^p(M)$.
For now the formulas have been proved for smooth sections. Under
approximation in $W^{r,p}(\mathbb R^m)$, the right hand sides stay bounded in $L^p(\mathbb R^m)$
by the preceding estimates. Since $1<p<\infty$, we can extract weak
limits for the finitely many derivatives considered. Uniform convergence
of the sections identifies the order zero limit with
\eqref{eq:derivadas-transicion-sobolev-explicitas}; integration against
test functions identifies the other limits with its distributional
derivatives. This proves membership in $W^{r,p}(\mathbb R^m)$ and the same bound for
Sobolev sections. We next prove the strong convergence needed for continuity.

We check continuity of $A_a$ in operator norm. For terms whose
coefficients have geometric order at most $k$, the lemma makes those
coefficients uniformly Lipschitz in $v$. The uniform embedding and the
product estimates allow passage to the limit when
$\mathbf X_n\to\mathbf X$ in $W^{r,p}(M,f^*TN)$. Only coefficients of order
$k+1$ remain; these can occur for $a=q$.

Every section $\mathbf U\in W^{r,p}(M,\mathbf F)$ of the bundles in question tends
uniformly to zero outside compact sets. Indeed, $\nabla\mathbf U$ is
bounded, so $|\mathbf U|$ is Lipschitz. If some $\varepsilon>0$ and
points escaping every compact set satisfied
$|\mathbf U(x_j)|\geq\varepsilon$, this bound would yield balls of
one positive radius on which $|\mathbf U|\geq\varepsilon/2$.
Completeness allows a subsequence of disjoint balls, and bounded
geometry gives a common positive lower bound for their volumes.
Then $\displaystyle \int_M|\mathbf U|^p\,d\lambda_{\mathbf g}=\infty$, a
contradiction. Applying this to $\mathbf X_0$ and $\mathbf Y_0$ also
shows that $d_{\mathbf h}(f(x),g(x))$ tends uniformly to zero outside
compact sets.

Let $\mathbf X_n\to\mathbf X$ in $W^{r,p}(M,f^*TN)$ and fix $\varepsilon>0$.
Outside a compact set $K\subset M$, both $|\mathbf X|$ and
$d_{\mathbf h}(f,g)$ are less than $\varepsilon$. For large $n$,
also $\|\mathbf X_n-\mathbf X\|_\infty<\varepsilon$.
By \eqref{eq:modulo-jets-comparacion-diagonal}, the difference of each
coefficient of order $k+1$, evaluated at $\mathbf X_n$ and $\mathbf X$,
is bounded on $M\setminus K$ by $C\varepsilon$. On $K$ it tends
uniformly to zero by smoothness on a compact domain. Taking first
$n\to\infty$ and then $\varepsilon\to0^{+}$ proves uniform
convergence of the coefficient on the whole base. Differences of the
remaining factors in \eqref{eq:terminos-transicion-indices-sobolev}
are estimated by H\"older's inequality as before. Applied also to smooth
approximations, this argument identifies the weak derivative formulas
and proves their strong convergence. Therefore
\[
 \|A_a(\mathbf X_n)-A_a(\mathbf X)\|_{\mathcal L^a(W^{r,p}(M,f^*TN);W^{r,p}(M,g^*TN))}
 \longrightarrow0,\qquad 1\leq a\leq q.
\]
This uses $p<\infty$ and the proximity of centers in an overlap;
it does not require uniform continuity of $\nabla^kR$ on all of $N$.

Since $\Psi(\,\cdot\,,\mathbf X_0)=\mathbf Y_0\in W^{r,p}(M,g^*TN)$,
the pointwise Taylor identity and continuity of $A_1$ give
\[
 \mathcal T(\mathbf X)-\mathbf Y_0
 =\int_0^1 A_1(\mathbf X_0+t(\mathbf X-\mathbf X_0))
                   [\mathbf X-\mathbf X_0]\,dt
 \quad\text{in }W^{r,p}(M,g^*TN).
\]
This Bochner integral exists. The identity proves membership in the
stated space and continuity of $\mathcal T$. For higher derivatives,
applying the same identity to $D_v^{a-1}\Psi$ gives
\[
 A_{a-1}(\mathbf X+\mathbf H)-A_{a-1}(\mathbf X)
 =\int_0^1 A_a(\mathbf X+t\mathbf H)
                   [\mathbf H,\,\cdot\,]\,dt,
 \qquad 2\leq a\leq q.
\]
After subtracting the integrand at $t=0$, the remainder norm is at most
$\|\mathbf H\|_{W^{r,p}(M,f^*TN)}$ times
$\displaystyle \sup_{t\in[0,1]}\|A_a(\mathbf X+t\mathbf H)-A_a(\mathbf X)\|$.
Given $\varepsilon>0$, operator-norm continuity provides
$\delta=\delta(\mathbf X,a,\varepsilon)>0$ such that the supremum is less
than $\varepsilon$ whenever $\|\mathbf H\|_{W^{r,p}(M,f^*TN)}<\delta$ and the segment
remains in the chart. Thus the remainder norm is at most
$\varepsilon\|\mathbf H\|_{W^{r,p}(M,f^*TN)}$. The same argument with $A_1$ proves
$D\mathcal T=A_1$, and induction gives
$D^a\mathcal T=A_a$ for $1\leq a\leq q$.

Shrink the neighborhood of $\mathbf X_0$ so that
$\|\mathcal T(\mathbf X)-\mathbf Y_0\|_\infty
<\rho-\|\mathbf Y_0\|_\infty$. This positive margin proves openness
of the overlap. Interchanging $f$ and $g$ produces an inverse transition
of the same class. Point evaluations are continuous; distinct mappings
are separated by disjoint neighborhoods of their values at some point.
The atlas topology is therefore Hausdorff. If geometry is bounded in
all orders, the argument applies to every $q\in\mathbb N$, and the
atlas is smooth.

The intrinsic identification of sections along a Sobolev mapping with
the model at a smooth center follows from transport in
Proposition~\ref{prop:transporte-homotopia-conexion-sobolev}.
Its proof uses the geometric estimates above and the module lemma,
without presupposing the atlas.
Exponential coordinates centered at $h$ also give topological
coordinates. Indeed, let $P\colon f^*TN\to h^*TN$ be that transport.
For a small section $\mathbf Z$ in the model at $f$, the coordinate
change and its inverse are, pointwise,
\[
 \mathbf Z\longmapsto
 P^{-1}\log_h\!\left(\operatorname{Exp}_f(\mathbf X_0+\mathbf Z)\right),
 \qquad
 \mathbf V\longmapsto
 \log_f\!\left(\operatorname{Exp}_h(P\mathbf V)\right)-\mathbf X_0.
\]
Both expressions vanish at the zero section. In trivializations at the
smooth center, they are superpositions in the two sections
$\mathbf X_0,\mathbf Z$, or $\mathbf X_0,\mathbf V$.
The integral difference formula and mixed bounds through order $r+1$
give membership in $W^{r,p}(M,f^*TN)$ and continuity by the same product
estimate. Only $k\geq r$ is needed here.
Their compositions are the pointwise identity, so they are inverse
homeomorphisms on sufficiently small neighborhoods. This justifies
the later use of $\operatorname{Exp}_h$ without treating $h^*TN$ as
a smooth bundle.

At a smooth center,
$D\Phi_f(0)\mathbf H=\mathbf H$ pointwise. Differentiating the inverse
chart identities gives the bounded linear isomorphisms in the statement.
\end{proof}

Smooth centers are locally dense: approximate a section
$\mathbf X\in\mathcal V_f$ by smooth compactly supported sections in
$W^{r,p}(M,f^*TN)$. Uniform embedding keeps the approximations inside
$\mathcal V_f$, and their images under $\operatorname{Exp}_f$ converge to the
initial map. Chart radii and transition constants
may vary with the center. This local assertion does not assign a
constant-radius uniformity to all maps with bounded derivatives.

\subsubsection{Parallel transport and the connection of a Sobolev map}
\label{subsubsec:transporte-conexion-mapeo-sobolev}

Identifying section spaces in two charts must not be confused with
an arbitrary change of connection. The comparison uses parallel
transport along the exponential homotopy. We make this mechanism
precise under the same finite-order hypotheses as the atlas, for
bases and targets without boundary.

\begin{proposition}[The connection difference induced by a homotopy]
\label{prop:transporte-homotopia-conexion-sobolev}
Let $M,N$ be without boundary, and of bounded geometry of
order $k$, and let $1<p<\infty$ and $k\geq r>m/p+1$, with
$k,r\in\mathbb N$. Fix
$f\in C_b^{\infty,r}(M,N)$ and denote the connection on
$\mathbf E=f^*TN$ by $\nabla^f$. Let $\mathbf X\in W^{r,p}(M,\mathbf E)$
have uniform norm smaller than a sufficiently small common normal radius.
Define
\[
 c_t(x)=\exp_{f(x)}(t\mathbf X(x)),\qquad h=c_1,
 \qquad P_t(x)\colon T_{f(x)}N\longrightarrow T_{c_t(x)}N
\]
by radial parallel transport, $0\leq t\leq1$.
The transported connection $P_1^{-1}\nabla^hP_1$ has the form
$\nabla^f+\mathbf a_{\mathbf X}$, where
\begin{equation}\label{eq:variacion-conexion-transporte-curvatura}
 \mathbf a_{\mathbf X}(\mathbf Y)\mathbf u
 =\int_0^1 P_t^{-1}R^N\bigl(\partial_t c_t,dc_t\mathbf Y\bigr)
                         P_t\mathbf u\,dt.
\end{equation}
Moreover,
\[
 \mathbf a_{\mathbf X}\in
 W^{r-1,p}(M,T^*M\otimes\operatorname{End}\mathbf E;\nabla^f),
\]
and $\mathbf X\mapsto\mathbf a_{\mathbf X}$ is continuous in these
norms. Its bounds on bounded sets may depend on the derivative bounds
for $f$ and on the geometric data. The transport $P_1$ identifies
$W^{r,p}(M,f^*TN)$ bicontinuously with the domain of sections along $h$
whose weak covariant derivatives up to order $r$ belong to $L^p(M,(T^*M)^{\otimes j}\otimes h^*TN)$.
\end{proposition}

\begin{proof}
First assume that $\mathbf X$ is smooth. Extend a base vector field
$\mathbf Y$ independently of $t$, and transport a section $\mathbf u$
by $\mathbf u_t=P_t\mathbf u$.
Since $\nabla_t\mathbf u_t=0$ and $[\partial_t,\mathbf Y]=0$, the
curvature identity gives
\[
 \nabla_t\nabla_{\mathbf Y}\mathbf u_t
   =R^N(\partial_t c_t,dc_t\mathbf Y)\mathbf u_t.
\]
Differentiating $P_t^{-1}\nabla_{\mathbf Y}\mathbf u_t$ amounts to
transporting this covariant derivative back to the fiber at $f(x)$.
At $t=0$ the result is $\nabla^f_{\mathbf Y}\mathbf u$.
Integration from zero to one gives
\eqref{eq:variacion-conexion-transporte-curvatura}, with the sign fixed
by our curvature convention.

We prove its Sobolev regularity. Put $a=r-1>m/p$.
Use uniform base charts small enough that $f$ takes values in normal
target charts of a common size. This is possible since
$\|df\|_\infty<\infty$. The bounds on derivatives of $f$ up to order $r$
control the induced connection coefficients to the orders in use.
Equivalence of coordinate and covariant norms follows from the same
derivative expansion as in
Proposition~\ref{prop:covariantes-localizacion-orden-entero}.

In these trivializations, the coefficients of $d\exp$, transport,
and transported curvature have the form $F(t,x,\mathbf X(x))$.
We count the orders needed. To control a superposition difference in
$W^{a,p}(\mathbb R^m)$, we use mixed derivatives of $F$ through order $a+1=r$.
For $d\exp$, these require derivatives of the exponential map through
order $r+1$; by Lemma~\ref{lem:jets-mixtos-exponencial-orden-finito},
bounds on $\nabla^jR^N$ for $0\leq j\leq r$ suffice.
The same orders suffice for transport. Differentiating transported
curvature $r$ times does not exceed $\nabla^rR^N$ either.
Derivatives of the centers only reach order $r$.
Thus $k\geq r$ supplies all the bounds, uniformly for $t\in[0,1]$.

More explicitly, write
\[
 F(t,x,\mathbf X(x))-F(t,x,\mathbf Z(x))
 =\int_0^1 D_vF(t,x,\mathbf Z(x)+s(\mathbf X(x)-\mathbf Z(x)))
                   [\mathbf X(x)-\mathbf Z(x)]\,ds.
\]
After differentiating through order $a$ in $x$, every term contains a
derivative of $\mathbf X-\mathbf Z$. The product estimate used in
\eqref{eq:cota-operador-transicion-sobolev}, now with $a>m/p$, gives
\[
 \|F(t,\,\cdot\,,\mathbf X)-F(t,\,\cdot\,,\mathbf Z)\|_{W^{a,p}(\mathbb R^m)}
 \leq C_B\|\mathbf X-\mathbf Z\|_{W^{a,p}(\mathbb R^m)}
\]
when the $W^{a,p}(\mathbb R^m)$ norms of $\mathbf X,\mathbf Z$ are less than $B$
and their values stay in the fixed normal domain.
The constant is independent of $t$. Setting $\mathbf Z=0$ gives
\[
 F(t,x,\mathbf X(x))=F(t,x,0)+F_{\mathbf X}(t,x),
 \qquad F_{\mathbf X}(t,\cdot)\in W^{a,p}(\mathbb R^m).
\]
Here $F(t,\cdot,0)$ has uniformly bounded derivatives up to order $a$,
whereas $F_{\mathbf X}$ depends continuously on $\mathbf X$ in
$W^{a,p}(\mathbb R^m)$, uniformly in $t$. The subtraction is essential: Taylor is
applied to the difference so that every term contains at least one
factor $\mathbf X$ or one of its derivatives. No constant is integrated
over a base of infinite volume. Passage from the charts to the global
norm uses bounded multiplicity, as in
Lemma~\ref{lem:modulo-sobolev-uniforme-entero}.

Each preceding coefficient is thus a bounded multiplier on $W^{a,p}(\mathbb R^m)$:
the part $F(t,\cdot,0)$ is handled by Leibniz and the part
$F_{\mathbf X}$ by the module lemma with $b=a$.
Moreover, $P_t^{-1}\partial_t c_t=\mathbf X$, since the radial velocity
is parallel. The horizontal and vertical decomposition of $d\exp$ gives
\[
 P_t^{-1}dc_t=B_t\,df+tA_t\nabla^f\mathbf X,
\]
where $A_t,B_t$ belong to the coefficient class just described.
After transporting $R^N$ as well, the integrand in
\eqref{eq:variacion-conexion-transporte-curvatura} is a sum of
contractions of multiplier coefficients with
$\mathbf X\otimes df$ and
$\mathbf X\otimes\nabla^f\mathbf X$.
The first product is in $W^{a,p}(M,T^*M\otimes\operatorname{End}\mathbf E)$ because derivatives of $df$ up to
order $a$ are bounded; the second is in $W^{a,p}(M,T^*M\otimes\operatorname{End}\mathbf E)$ by the module lemma.
These estimates are uniform in $t$ and continuous in
$\mathbf X\in W^{r,p}(M,f^*TN)$. The integral is therefore the integral of a
continuous curve with values in the Banach space $W^{a,p}(M,T^*M\otimes\operatorname{End}\mathbf E)$.
This proves the stated regularity and continuity.

For a merely Sobolev $\mathbf X$, choose smooth approximations in
$W^{r,p}(M,f^*TN)$. The uniform embedding into $C^1$ ensures convergence of the
homotopies, transports and induced connections in their first-order
coordinate expressions. The preceding estimates give convergence of
their connection differences in $W^{r-1,p}(M,T^*M\otimes\operatorname{End}\mathbf E)$. The limit agrees with the
difference of the $C^0$ connections induced by the $C^1$ map $h$, and
retains \eqref{eq:variacion-conexion-transporte-curvatura}.

Finally, transport is an isometry in each fiber. Under identification
by $P_1$, covariant derivatives along $h$ are those of
$\nabla^f+\mathbf a_{\mathbf X}$.
Corollary~\ref{cor:conexion-sobolev-dominios-debiles} compares its
domain with $W^{r,p}(M,f^*TN;\nabla^f)$ and its norms in both directions.
The proof of that corollary only needs the module lemma for the reference,
which is available in the preceding trivializations; it does not demand
bounded geometry of a connection that may now only have Sobolev
coefficients. This gives the required bicontinuous identification.
\end{proof}

This comparison supplements the local atlas construction. In particular,
it does not identify arbitrary pairs of connections on a noncompact base
without hypotheses. It explains why the specific geometric transport
relating chart centers provides the required estimates. The preceding
count shows that curvature bounds through order $r$ suffice to compare
the weak domains of order $r$; bounded geometry of all orders has not
been assumed.

\section{Spaces of metrics, connections, and gauge theory}
\label{sec:espacios-metricas-conexiones-gauge}

\subsection{Riemannian metrics}

Let $M$ be a smooth compact manifold of dimension $m$ with or without boundary and let
$s>\frac{m}{2}$. Define
\[
\operatorname{Met}^s(M):=
\{\mathbf{g}\in H^s(M,S^2T^*M)\mid
\mathbf{g}_x\text{ is positive definite for every }x\in M\}.
\]

\begin{proposition}[Openness of positivity]
\label{prop:metrica-positividad-abierta}
Let $M$ be a smooth compact manifold of dimension $m$ with or without boundary and let
$s>\frac{m}{2}$. Then $\operatorname{Met}^s(M)$ is an open subset of the Hilbert space
$H^s(M,S^2T^*M)$. In particular, it is a Hilbert manifold and
$T_{\mathbf{g}}\operatorname{Met}^s(M)=H^s(M,S^2T^*M)$.
\end{proposition}

\begin{proof}
By compactness, the smallest eigenvalue of $\mathbf{g}$ relative to an auxiliary metric
$\mathbf{g}_0$ has a minimum $c_{\mathbf{g}}>0$. The embedding
\[
H^s(M,S^2T^*M)\hookrightarrow C^0(M,S^2T^*M;\mathbf{g}_0)
\]
gives $C>0$ such that
\[
\|\mathbf{h}\|_{C^0(M,S^2T^*M;\mathbf{g}_0)}
\leq C\|\mathbf{h}\|_{H^s(M,S^2T^*M)}.
\]
If
$\|\mathbf{h}\|_{H^s(M,S^2T^*M)}<\frac{c_{\mathbf{g}}}{2C}$, then
$(\mathbf{g}+\mathbf{h})(v,v)\geq(\frac{c_{\mathbf{g}}}{2})|v|_{\mathbf{g}_0}^2$
for every $v\in T_xM$. Thus the corresponding ball
lies in $\operatorname{Met}^s(M)$. Since this set is open, its tangent
space is the model space.
\end{proof}

Example~\ref{ej:metricas-conformes-banach} fixed a metric
$\mathbf g_0$ and constructed the family $a\mathbf g_0$ with
strictly positive $a\in C^k(M)$. Its tangent directions are
conformal to the metric at the base point. In $\operatorname{Met}^s(M)$,
regularity is measured in the $H^s$ norm and tangent directions
range over all symmetric tensors in that space, as the preceding
proposition shows.

On a noncompact manifold, the smallest eigenvalue may tend to zero at
infinity; positivity is not open in an absolute Sobolev topology. For
example, on $\mathbb R$ the metric $\mathbf{g}=e^{-2x^2}(\mathbf{d}x)^2$ is positive, but there exist
cutoff functions $\chi_j$ translated toward $j$, with uniformly bounded
derivatives, such that
$\mathbf{h}_j=-2e^{-2x^2}\chi_jdx^2$ tends to zero in every $H^s(\mathbb R,S^2T^*\mathbb R)$ while
$\mathbf{g}+\mathbf{h}_j$ is negative where $\chi_j=1$. This difficulty disappears when working
within a uniform equivalence class.

Two metrics $\mathbf{g}$ and $\mathbf{g}'$ are \textit{quasi-isometric} if there exists $C\geq1$
such that
\begin{equation}
\label{eq:quasi-isometria-metricas-puntual}
 C^{-1}\mathbf{g}\leq \mathbf{g}'\leq C\mathbf{g}
\end{equation}
in the sense of Definition~\ref{def:orden-formas-metricas};
that is, a common constant $C$ satisfies
\[
 C^{-1}\mathbf g_x(v,v)\leq\mathbf g'_x(v,v)
 \leq C\mathbf g_x(v,v)
 \qquad(x\in M,\ v\in T_xM).
\]
Equivalently, the tangent norms satisfy
\[
 C^{-1/2}|v|_{\mathbf g}\leq|v|_{\mathbf g'}
 \leq C^{1/2}|v|_{\mathbf g}
 \qquad(x\in M,\ v\in T_xM).
\]
This condition controls the values of the metrics on all
fibers with common constants; control of the connections is
added in the following functional.
If $\nabla^{\mathbf{g}}$ and $\nabla^{\mathbf{g}'}$ are their Levi--Civita
connections, define the comparison functional
\begin{equation}
\label{eq:norma-bq-diferencia-metricas}
 \delta_{b,q;\mathbf{g}}(\mathbf{g},\mathbf{g}')
 :=\|\mathbf g-\mathbf g'\|_{\infty;\mathbf g}+
 \sum_{j=0}^{q-1}
 \|(\nabla^{\mathbf{g}})^j(\nabla^{\mathbf{g}'}-\nabla^{\mathbf{g}})\|_{\infty;\mathbf g}.
\end{equation}
The infinity norms use the metrics on the indicated bundles;
when confusion could arise, the metric is added after a semicolon.
The difference of connections is a tensor of type $(1,2)$, so all
terms are intrinsically defined.

For $k\in\mathbb N_0$, denote by
$\operatorname{Met}(I,B_k)(M)$ the set of smooth Riemannian
metrics $\mathbf{g}$ satisfying conditions $(I)$ and $(B_k)$ in
Definition~\ref{def:geometria-acotada-orden-k}; that is, there exist constants
$i_{\mathbf{g}},K_0(\mathbf{g}),\ldots,K_k(\mathbf{g})>0$ such that
\[
 \operatorname{inj}(M,\mathbf{g})\geq i_{\mathbf{g}},
 \qquad
 \|(\nabla^{\mathbf{g}})^j\mathbf{R}_{\mathbf{g}}\|_{\infty;\mathbf g}\leq K_j(\mathbf{g})
 \quad(0\leq j\leq k).
\]
The first condition implies that $(M,\mathbf{g})$ is complete. In the following
Sobolev construction, both metrics of each pair belong to this set.

If $\mathbf{T}$ is a tensor field of fixed type and $j\in\mathbb N_0$, we use
$\displaystyle \|\mathbf{T}\|_{b,j;\mathbf{g}}:=\displaystyle\sum_{a=0}^j\|(\nabla^{\mathbf{g}})^a\mathbf{T}\|_{\infty;\mathbf g}$. For
$1<p<\infty$, we also set
\[
 \|\mathbf{T}\|_{p,j;\mathbf{g}}
 :=\left(
 \sum_{a=0}^j\int_M|(\nabla^{\mathbf{g}})^a\mathbf{T}|_{\mathbf{g}}^p\,d\lambda_{\mathbf{g}}
 \right)^{\frac{1}{p}}.
\]
Both notations are used only when the right-hand side is
finite.

\begin{theorem}[Uniform comparison of metrics]
\label{prop:uniformidad-metricas-eichhorn}
Let $M$ (of dimension $m$) be smooth and without boundary. For $q\in\mathbb N$, the sets
\[
 \mathcal U^{b,q}_\varepsilon
 :=\left\{(\mathbf g,\mathbf g')\mid
 \mathbf g\text{ and }\mathbf g'\text{ are quasi-isometric and }
 \delta_{b,q;\mathbf g}(\mathbf g,\mathbf g')<\varepsilon\right\}
\]
generate a separated metrizable uniformity on the smooth complete
metrics. In the Sobolev regime, fix a reference metric
$\mathbf g_0\in\operatorname{Met}(I,B_k)(M)$, $1<p<\infty$, and
$k\geq r>m/p+1$. The intrinsic comparison is written as
\begin{equation}
\label{eq:norma-pr-diferencia-metricas}
 \delta_{p,r;\mathbf g}(\mathbf g,\mathbf g')^p
 :=\int_M\left(|\mathbf g-\mathbf g'|_{\mathbf g}^p+
 \sum_{j=0}^{r-1}|(\nabla^{\mathbf g})^j(\nabla^{\mathbf g'}-\nabla^{\mathbf g})|_{\mathbf g}^p\right)
 \,d\lambda_{\mathbf g}.
\end{equation}
On the affine space $\mathbf g_0+W^{r,p}(M,S^2T^*M;\mathbf g_0)$,
we use the distance
\[
 d_{\mathbf g_0}(\mathbf g,\mathbf g')
       :=\|\mathbf g'-\mathbf g\|_{p,r;\mathbf g_0}.
\]
On every set of smooth metrics satisfying, for fixed numbers
$\lambda>0$ and $R>0$,
\[
 \mathbf g\geq\lambda\mathbf g_0,
 \qquad \|\mathbf g-\mathbf g_0\|_{p,r;\mathbf g_0}\leq R,
\]
there exist $c,C>0$ common to all pairs such that
\[
 c\,d_{\mathbf g_0}(\mathbf g,\mathbf g')
 \leq\delta_{p,r;\mathbf g}(\mathbf g,\mathbf g')
 \leq C\,d_{\mathbf g_0}(\mathbf g,\mathbf g').
\]
The constants depend on $\lambda,R,p,r$ and the data of $\mathbf g_0$.
\end{theorem}

\begin{proof}
If $\|\mathbf g'-\mathbf g\|_{\infty;\mathbf g}<1/2$, for
every $x\in M$ and $v\in T_xM$, Cauchy--Schwarz
in the tensor inner product gives
\[
 |(\mathbf g'-\mathbf g)_x(v,v)|
 \leq|\mathbf g'-\mathbf g|_{\mathbf g}(x)|v|_{\mathbf g}^2
 \leq\frac12\mathbf g_x(v,v).
\]
Adding $\mathbf g_x(v,v)$ gives
\[
 \frac12\mathbf g_x(v,v)\leq\mathbf g'_x(v,v)
 \leq\frac32\mathbf g_x(v,v).
\]
By Proposition~\ref{prop:comparacion-metricas-duales-tensores-volumen},
the pointwise tensor norms are comparable with constants
depending only on the tensor type. The change-of-connection formula
\eqref{eq:expansion-cambio-conexion-tensor-metricas}, proved
below, expresses derivatives in the new connection as finite sums
of products of derivatives of the connection difference. Applied to
$\mathbf g'-\mathbf g$ and $\nabla^{\mathbf g'}-\nabla^{\mathbf g}$,
it gives polynomials with no constant term controlling inversion.
For three metrics, add their connection differences and apply the
same formula when changing connection in the second summand. This gives
corresponding composition polynomials independent of the
metrics. This verifies the inversion and composition axioms; the other
axioms follow from the diagonal and order-zero control. Rational
radii and metrization by chains complete the argument, also
formulated in \cite[Proposition~2.1]{EichhornMetrics}.

For Sobolev norms, the embedding gives
$\|\mathbf g-\mathbf g_0\|_{\infty;\mathbf g_0}\leq C_SR$.
The same evaluation on $(v,v)$ gives
$\mathbf g_x(v,v)\leq(1+C_SR)(\mathbf g_0)_x(v,v)$
for every $x$ and $v$. Together with the lower bound in the statement and
Proposition~\ref{prop:comparacion-metricas-duales-tensores-volumen},
this implies
\[
 (1+C_SR)^{-1}\mathbf g_0^{-1}\leq\mathbf g^{-1}
 \leq\lambda^{-1}\mathbf g_0^{-1},
 \qquad
 \lambda^{m/2}
 \leq\frac{d\lambda_{\mathbf g}}{d\lambda_{\mathbf g_0}}
 \leq(1+C_SR)^{m/2}.
\]
Comparison of the inverse metrics is evaluated on covectors,
and the ratio of measures is the positive function described in
that proposition. These bounds are common to all metrics
in the set under consideration.
The Koszul formula \eqref{eq:diferencia-levi-civita-metricas} expresses the
connection difference in terms of one metric's inverse and the first
derivative of the other. Differentiating the identity
$\mathbf g^{-1}\mathbf g=I$ gives
$\nabla\mathbf g^{-1}=-\mathbf g^{-1}(\nabla\mathbf g)\mathbf g^{-1}$;
at each subsequent order, Leibniz's rule differentiates one of the existing factors.
Thus, through order $r-1$, derivatives of the inverse and the connection
difference are controlled by $\lambda,R$ and the product estimates.
In particular, $\nabla^{\mathbf g}-\nabla^{\mathbf g_0}$ is bounded in
$W^{r-1,p}(M,T^*M\otimes\operatorname{End}TM)$ in terms of these data.

Set $\mathbf h=\mathbf g'-\mathbf g$ and
$\mathbf D=\nabla^{\mathbf g'}-\nabla^{\mathbf g}$. Koszul's formula and the product estimates
give $\|\mathbf D\|_{p,r-1;\mathbf g_0}\leq C\|\mathbf h\|_{p,r;\mathbf g_0}$.
For the reverse inequality, use the metric compatibility identity
\[
 (\nabla^{\mathbf g}_{\mathbf X}\mathbf h)(\mathbf Y,\mathbf Z)
 =\mathbf g'(\mathbf D(\mathbf X,\mathbf Y),\mathbf Z)
   +\mathbf g'(\mathbf Y,\mathbf D(\mathbf X,\mathbf Z)).
\]
Starting at order zero for $\mathbf h$, differentiate this identity through
order $r-1$. In the step from order $j$ to $j+1$, terms in which
$\mathbf g'=\mathbf g+\mathbf h$ is differentiated contain derivatives of
$\mathbf h$ of order at most $j$, already controlled; the remaining terms contain
a derivative of $\mathbf D$ of order at most $j$. The module estimate
for $W^{r-1,p}(M,\mathbf F)$ in each tensor bundle $\mathbf F\to M$ involved, valid because $r-1>m/p$, controls all these products.
The change-of-connection formula allows each norm to be expressed using
$\mathbf g_0$, and the volume ratios are bounded. We obtain
$\|\mathbf h\|_{p,r;\mathbf g_0}\leq
C(\|\mathbf h\|_{p,0;\mathbf g}+\|\mathbf D\|_{p,r-1;\mathbf g})$,
which is the remaining inequality.
\end{proof}

The fixed reference is essential in the Sobolev assertion. On $S^1$, with
$\theta\in\mathbb R/(2\pi\mathbb Z)$, take
$\mathbf g_a=a^2d\theta^2$ and $\mathbf h_a=a^6d\theta^2$, for $0<a<1$.
Their connections agree, and both metrics are complete and of bounded
geometry, though their data are not common as $a$ varies. For $p=2$ and any
integer order $r\geq1$,
\[
 \delta_{2,r;\mathbf g_a}(\mathbf g_a,\mathbf h_a)
 =\sqrt{2\pi}\,a^{1/2}(1-a^4),
 \qquad
 \delta_{2,r;\mathbf h_a}(\mathbf h_a,\mathbf g_a)
 =\sqrt{2\pi}\,a^{-5/2}(1-a^4).
\]
The first side tends to zero and the second to infinity. Thus the condition
$\delta_{p,r;\mathbf g}<\varepsilon$ on all metrics without common
data does not satisfy the inversion axiom. The distance
$d_{\mathbf g_0}$ is symmetric and satisfies the triangle inequality; the
preceding comparison allows the intrinsic expressions to be used in each
neighborhood with fixed bounds.

\begin{proposition}[Comparison of tensor norms in a component]
\label{prop:comparacion-normas-metricas-componente}
Let $M$ be a smooth manifold without boundary, let $\mathbf{g}$ and $\mathbf{g}'$ be
smooth quasi-isometric metrics, and let $\mathbf{T}$ be a tensor field
of fixed type.
\begin{enumerate}[label=(\roman*)]
\item If $\delta_{b,q;\mathbf g}(\mathbf g,\mathbf g')<\infty$, then, for
each $0\leq j\leq q$, there exist constants $c_{b,j},C_{b,j}>0$ such that
\[
 c_{b,j}\|\mathbf{T}\|_{b,j;\mathbf{g}}
 \leq\|\mathbf{T}\|_{b,j;\mathbf{g}'}
 \leq C_{b,j}\|\mathbf{T}\|_{b,j;\mathbf{g}}.
\]
\item If $\mathbf g,\mathbf g'$ belong to the affine space based at
$\mathbf g_0$ in the preceding theorem and have a positive lower bound
relative to $\mathbf g_0$, then, for
each $0\leq j\leq r$, there exist constants $c_{p,j},C_{p,j}>0$ such that
\[
 c_{p,j}\|\mathbf{T}\|_{p,j;\mathbf{g}}
 \leq\|\mathbf{T}\|_{p,j;\mathbf{g}'}
 \leq C_{p,j}\|\mathbf{T}\|_{p,j;\mathbf{g}}.
\]
\end{enumerate}
The constants depend on the tensor type, a quasi-isometry
constant between $\mathbf{g}$ and $\mathbf{g}'$, and the norms of
$(\nabla^{\mathbf{g}})^a(\nabla^{\mathbf{g}'}-\nabla^{\mathbf{g}})$ through the required order, but not on
$\mathbf{T}$.
\end{proposition}

\begin{proof}
Set $\mathbf{D}:=\nabla^{\mathbf{g}'}-\nabla^{\mathbf{g}}$. The Koszul identity gives
\begin{equation}
\label{eq:diferencia-levi-civita-metricas}
 2\mathbf{g}'\bigl(\mathbf{D}(\mathbf{X},\mathbf{Y}),\mathbf{Z}\bigr)
 =(\nabla^{\mathbf{g}}_{\mathbf{X}}\mathbf{g}')(\mathbf{Y},\mathbf{Z})
 +(\nabla^{\mathbf{g}}_{\mathbf{Y}}\mathbf{g}')(\mathbf{X},\mathbf{Z})
 -(\nabla^{\mathbf{g}}_{\mathbf{Z}}\mathbf{g}')(\mathbf{X},\mathbf{Y}).
\end{equation}
If $C^{-1}\mathbf{g}\leq \mathbf{g}'\leq C\mathbf{g}$,
Proposition~\ref{prop:comparacion-metricas-duales-tensores-volumen},
applied with constants $C^{-1}$ and $C$, gives, for a tensor of
type $(a,b)$,
\begin{equation}
\label{eq:comparacion-puntual-tensores-metricas}
 C^{-(a+b)/2}|\mathbf{T}|_{\mathbf{g}}\leq|\mathbf{T}|_{\mathbf{g}'}
 \leq C^{(a+b)/2}|\mathbf{T}|_{\mathbf{g}},
\end{equation},
and the ratio of volume elements satisfies
\[
 C^{-m/2}\leq\rho(x):=
 \frac{\sqrt{\det(\mathbf g'_x)}}{\sqrt{\det(\mathbf g_x)}}
 \leq C^{m/2}.
\]
The determinants are computed in a common chart, and the ratio
is independent of it. For $1\leq p<\infty$, we have
\[
 \|\mathbf T\|_{L^p(M,\mathbf g')}^p
 =\int_M|\mathbf T|_{\mathbf g'}^p\rho\,d\lambda_{\mathbf g}.
\]
Raising \eqref{eq:comparacion-puntual-tensores-metricas} to the
power $p$, multiplying by the bounds on $\rho$, and integrating
gives
\[
 C^{-(a+b)/2-m/(2p)}\|\mathbf T\|_{L^p(M,\mathbf g)}
 \leq\|\mathbf T\|_{L^p(M,\mathbf g')}
 \leq C^{(a+b)/2+m/(2p)}\|\mathbf T\|_{L^p(M,\mathbf g)}.
\]
For the uniform norm, simply take suprema in the pointwise
comparison. These estimates compare the same tensor $\mathbf T$;
the next step compares the derivatives obtained using the two
connections.

The local expansion of the covariant derivative in
Lemma~\ref{lem:local-expression-higher-order} and the Leibniz rule in
Proposition~\ref{prop: leibniz operador *} allow us to write, for
$j\in\mathbb N$,
\begin{equation}
\label{eq:expansion-cambio-conexion-tensor-metricas}
 (\nabla^{\mathbf{g}'})^j\mathbf{T}=(\nabla^{\mathbf{g}})^j\mathbf{T}+
 \sum_{\ell=1}^j
 \sum_{i_0+\cdots+i_\ell=j-\ell}
 C_{j,\boldsymbol i}
 (\nabla^{\mathbf{g}})^{i_1}\mathbf{D}*\cdots*(\nabla^{\mathbf{g}})^{i_\ell}\mathbf{D}*
 (\nabla^{\mathbf{g}})^{i_0}\mathbf{T}.
\end{equation}
In the inner sum, $\boldsymbol i=(i_0,\ldots,i_\ell)$ ranges over
$\mathbb N_0^{\ell+1}$ subject to the indicated restriction. The coefficients are
universal; each term with $*$ represents the finite
sum of contractions and permutations specified by the tensor type. The
endomorphisms appearing here are not assumed to commute. For $j=1$, we obtain
$\nabla^{\mathbf{g}'}\mathbf{T}=\nabla^{\mathbf{g}}\mathbf{T}+\mathbf{D}*\mathbf{T}$;
on differentiating again, the derivative may act
on any existing factor or introduce a new factor
$\mathbf{D}$. In the first case, one of the $i_a$, with
$a\in\{0,\ldots,\ell\}$, increases by one, and the index sum changes
from $j-\ell$ to $(j+1)-\ell$. In the second, $\ell$ increases by one and a
zero index is added; the sum remains equal to
$(j+1)-(\ell+1)$. These are all the terms of order $j+1$, and their orders
and types are those in the formula. This completes the induction step.

In the $b,q$ uniformity, supremum bounds on the derivatives of $\mathbf{D}$ and
\eqref{eq:comparacion-puntual-tensores-metricas} control
\eqref{eq:expansion-cambio-conexion-tensor-metricas} term by term. In the Sobolev
uniformity, $r>m/p+1$ implies that $W^{r-1,p}(M)$ is an algebra and its sections
act as multipliers on lower orders. The product and
module theorems, applied to each contraction, give the upper bound
in part (ii). Interchanging $\mathbf{g}$ and $\mathbf{g}'$ gives the reverse inequalities
in both cases.
\end{proof}

For fixed $\mathbf g_0\in\operatorname{Met}(I,B_k)(M)$, denote by
$\operatorname{Met}_{\mathbf g_0}^{p,r}(M)$ the set
\[
 \left\{\mathbf g_0+\mathbf h\ \middle|\
 \mathbf h\in W^{r,p}(M,S^2T^*M;\mathbf g_0),\quad
 \mathbf g_0+\mathbf h\geq c\mathbf g_0\text{ for some }c>0\right\}.
\]
This is the positive affine component considered when writing
$\operatorname{Met}^{p,r}(I,B_k)(M)$ with a fixed reference. The complete affine
space is $\mathbf g_0+W^{r,p}(M,S^2T^*M)$; uniform positivity selects
an open subset of it, which need not be complete for the affine distance.
In the uniform version, let $B^q_{\mathbf g_0}$ be the completion of the
smooth sections bounded through order $q$ in $\|\cdot\|_{b,q;\mathbf g_0}$.
Define $\operatorname{Met}^{b,q}_{\mathbf g_0}(M)$ by requiring
$\mathbf h\in B^q_{\mathbf g_0}$ and the same uniform positive inequality.

\begin{theorem}[Affine manifolds of metrics]
\label{teo:eichhorn-espacio-metricas}
Let $M$ (of dimension $m$) be smooth, and without boundary with respect to
$\mathbf g_0\in\operatorname{Met}(I,B_k)(M)$. If $1<p<\infty$ and
$k\geq r>m/p+1$, the space $\operatorname{Met}^{p,r}_{\mathbf g_0}(M)$ is
an open convex subset of $\mathbf g_0+W^{r,p}(M,S^2T^*M;\mathbf g_0)$.
It is a smooth Banach manifold, and a Hilbert manifold if $p=2$. For the affine chart
$\Phi_{\mathbf g}(\mathbf g')=\mathbf g'-\mathbf g$,
\[
 J_{\mathbf g}:=D\Phi_{\mathbf g}(\mathbf g):
 T_{\mathbf g}\operatorname{Met}^{p,r}_{\mathbf g_0}(M)
 \longrightarrow H^{r,p}(M,S^2T^*M;\mathbf g_0)
\]
is the linear bicontinuous identity. Here we use equivalence of the
$W^{r,p}(M,S^2T^*M)$ and $H^{r,p}(M,S^2T^*M)$ norms for integer orders from the bounded geometry
chapter. The uniform construction is an open convex subset modeled on
$B^q_{\mathbf g_0}$.
\end{theorem}

\begin{proof}
If $\mathbf g\geq c\mathbf g_0$, the uniform embedding gives
$\|\mathbf k\|_{\infty;\mathbf g_0}\leq C_S\|\mathbf k\|_{p,r;\mathbf g_0}$.
The condition $\|\mathbf k\|_{p,r;\mathbf g_0}<c/(2C_S)$ implies
$\mathbf g+\mathbf k\geq(c/2)\mathbf g_0$, proving openness. A
convex combination of two metrics with lower bounds $c_1,c_2>0$
has lower bound $(1-t)c_1+tc_2>0$ for $t\in[0,1]$; thus the
set is convex. Changes between affine charts are translations with
identity differential. The uniform assertion is proved in the same way, using the fact that
the $B^q_{\mathbf g_0}$ norm controls the infinity norm.

In the Sobolev component, pointwise positivity already implies a uniform
lower bound. To see this, $\mathbf h\in W^{r,p}(M,S^2T^*M)$ has bounded first
derivative and tends uniformly to zero at infinity. If the
latter assertion failed, there would exist $\varepsilon>0$ and points $x_j$
escaping every compact set, with $|\mathbf h(x_j)|\geq\varepsilon$. The bound on
$\nabla\mathbf h$ and parallel transport would give
$|\mathbf h|\geq\varepsilon/2$ on balls of a common positive radius around
these points. A disjoint family of balls can be extracted; their
volumes have a common lower bound by bounded geometry. Then
the sum of the integrals of $|\mathbf h|^p$ would be infinite, a contradiction.
Thus, outside a compact set,
$\mathbf g_0+\mathbf h\geq\frac12\mathbf g_0$. Inside it, continuity
and positivity give a positive minimum for the smallest eigenvalue.

Smooth compactly supported sections approximate $\mathbf h$ in
$W^{r,p}(M,S^2T^*M)$. Uniform convergence ensures that the approximations retain
a positive lower bound. The resulting smooth metrics agree with
$\mathbf g_0$ outside a compact set, are complete, and satisfy $(I)$ and
$(B_k)$, with constants that may depend on the approximation. This
relates the affine model to approximation by smooth metrics. It does not assert
that the limit retains curvature bounds of order $k$: that would
require control of additional metric derivatives.
\end{proof}

Affine charts are those used in studying spaces of metrics
on complete manifolds without boundary equipped with a bounded-geometry
reference metric; see \cite{EichhornMetrics}.
Proposition~\ref{prop:comparacion-normas-metricas-componente} explains why,
within each neighborhood with common bounds, one may use intrinsic norms
for another metric in the component. The fixed reference metric
makes explicit the topology used for completion and
openness of positivity.

On a compact base, or in a Hilbert Sobolev component of a
noncompact base, the Ebin metric is defined by
\begin{equation}
\label{eq:metrica-ebin-espacio-metricas}
G_{\mathbf{g}}(\mathbf{h},\mathbf{k})
:=\int_M\operatorname{tr}_{\mathbf{g}}(\mathbf{h}\circ\mathbf{k})\,d\lambda_{\mathbf{g}}.
\end{equation}
In the noncompact case, tangent vectors belong to $L^2(M,S^2T^*M)$ and the
metrics are uniformly equivalent to the reference; Hölder's inequality thus ensures
that the integral is finite. For a $W^{r,p}(M,S^2T^*M)$ component with $p\neq2$,
this integrability must additionally be required.
In an $H^s(M,S^2T^*M)$ topology with $s>0$, this metric controls only the $L^2(M,S^2T^*M)$ norm and
is therefore weak. Existence of a connection or an exponential map cannot
be deduced from the Koszul formula alone.

\subsection{Connections and gauge transformations}

Let $\mathbf{P}\longrightarrow M$ be a principal bundle with compact structure group
$G$ and Lie algebra $\mathfrak g$ equipped with a
$\operatorname{Ad}$--invariant inner product. Denote by $h_{\mathrm{ad}}$ the induced metric
on the adjoint bundle and by
$\operatorname{ad}\mathbf{P}:=\mathbf{P}\times_{\operatorname{Ad}}\mathfrak g$ the adjoint bundle.

\begin{proposition}[Affine structure of the space of connections]
\label{prop:espacio-conexiones-afin}
Let $M$ be a smooth manifold with or without boundary and let
$\mathbf{P}\longrightarrow M$ be a smooth principal bundle. For a fixed smooth
connection $A_0$, every smooth connection can be written
uniquely as $A_0+\mathbf{a}$, with
$\mathbf{a}\in\Omega^1(M,\operatorname{ad}\mathbf{P})$. Consequently, for every Sobolev
index $s$ for which the completion has been fixed,
\[
 \mathcal A^s(\mathbf{P})
 :=A_0+H^s(M,T^*M\otimes\operatorname{ad}\mathbf{P})
\]
is an affine space and
\[
 T_A\mathcal A^s(\mathbf{P})
 =H^s(M,T^*M\otimes\operatorname{ad}\mathbf{P}).
\]
\end{proposition}

\begin{proof}
The smooth assertion is
Theorem~\ref{teo:espacio-afin-conexiones-principales}: the difference of two
principal connection forms is horizontal and
$\operatorname{Ad}$--equivariant and therefore descends to a unique form
with values in $\operatorname{ad}\mathbf{P}$. Completing the model vector
space in $H^s(M,T^*M\otimes\operatorname{ad}\mathbf P)$ gives the second assertion.
\end{proof}

\begin{proposition}[Right gauge action and curvature]
\label{prop:accion-gauge-derecha-curvatura}
Let $M$ be a smooth manifold with or without boundary and let
$\mathbf{P}\longrightarrow M$ be a smooth principal bundle. The convention fixed in
\eqref{eq:accion-gauge-derecha-potencial} defines a right action:
\begin{equation}
\label{eq:accion-gauge-conexion}
 A\triangleleft \mathbf{u}
 :=\operatorname{Ad}_{\mathbf{u}^{-1}}A+\mathbf{u}^{-1}d\mathbf{u},
 \qquad
 (A\triangleleft \mathbf{u})\triangleleft \mathbf{v}=A\triangleleft(\mathbf{u}\mathbf{v}).
\end{equation}
In a matrix representation, the first formula is written as
$A\triangleleft \mathbf{u}=\mathbf{u}^{-1}A\mathbf{u}+\mathbf{u}^{-1}d\mathbf{u}$. Moreover,
\begin{align}
 \mathbf{F}_{A+\mathbf{a}}
 &=\mathbf{F}_A+d_A\mathbf{a}+\frac12[\mathbf{a}\wedge \mathbf{a}],
 \label{eq:curvatura-perturbacion-conexion}\\
 \mathbf{F}_{A\triangleleft \mathbf{u}}
 &=\operatorname{Ad}_{\mathbf{u}^{-1}}\mathbf{F}_A.
 \label{eq:curvatura-conexion-gauge}
\end{align}
For the Sobolev extension on a compact base, fix $s>m/2$;
on a noncompact base, impose the assumptions for the gauge component
constructed below. If
\[
 \Phi\colon\mathcal A^s(\mathbf{P})\times\mathcal G^{s+1}(\mathbf{P})
 \longrightarrow\mathcal A^s(\mathbf{P}),
 \qquad \Phi(A,\mathbf{u})=A\triangleleft \mathbf{u},
\]
and we trivialize $T_{\mathbf{u}}\mathcal G^{s+1}(\mathbf{P})$ by left translation,
writing a tangent vector as $\mathbf{u}\boldsymbol{\xi}$, then
\begin{equation}
\label{eq:diferencial-completo-accion-gauge-derecha}
 D\Phi_{(A,\mathbf{u})}(\mathbf{a},\mathbf{u}\boldsymbol{\xi})
 =\operatorname{Ad}_{\mathbf{u}^{-1}}\mathbf{a}+d_{A\triangleleft \mathbf{u}}\boldsymbol{\xi}.
\end{equation}
\end{proposition}

\begin{proof}
The first two identities are
Propositions~\ref{prop:formulas-gauge-conexiones-principales} and~\ref{teo:dA-cuadrado-bianchi-principal}, written in a local section. In
particular, the order of $\mathbf{u}\mathbf{v}$ in \eqref{eq:accion-gauge-conexion} is not an
implicit convention: it comes from the right action already established.

To compute the differential, take a curve $A(t)$ with
$A(0)=A$, $\dot A(0)=\mathbf{a}$, and the exact curve
$\mathbf{u}(t)=\mathbf{u}\exp_G(t\boldsymbol{\xi})$. By the action law,
\[
 A(t)\triangleleft \mathbf{u}(t)
 =\bigl(A(t)\triangleleft \mathbf{u}\bigr)\triangleleft\exp_G(t\boldsymbol{\xi}).
\]
The variation of the first factor is
$\operatorname{Ad}_{\mathbf{u}^{-1}}\mathbf{a}$, while the derivative of
$B\triangleleft\exp_G(t\boldsymbol{\xi})$ at $t=0$ is $d_B\boldsymbol{\xi}$. Taking
$B=A\triangleleft \mathbf{u}$ gives
\eqref{eq:diferencial-completo-accion-gauge-derecha}.
\end{proof}

\begin{proposition}[Coulomb decomposition in the compact case]
\label{prop:descomposicion-coulomb-compacta}
Let $(M,\mathbf{g})$ be a compact Riemannian manifold without boundary, let
\(\mathbf{P}\to M\) be a principal bundle with compact structure group, let
\(\mathbf{A}\) be a smooth connection, and let \(s\geq1\). Then
\[
 H^s\left(
 M,T^*M\otimes\operatorname{ad}\mathbf{P}
 \right)
 =
 \ker d_{\mathbf{A}}^*\oplus
 d_{\mathbf{A}}\left(
 H^{s+1}(M,\operatorname{ad}\mathbf{P})
 \cap(\ker d_{\mathbf{A}})^\perp
 \right).
\]
The two projections associated with this direct sum are continuous.
\end{proposition}

\begin{proof}
The operator
\[
 \Delta_{\mathbf{A}}^{(0)}:=d_{\mathbf{A}}^*d_{\mathbf{A}}
\]
is elliptic, nonnegative, and self-adjoint. Its resolvent maps
\(L^2(M,\operatorname{ad}\mathbf{P})\) continuously into
\(H^2(M,\operatorname{ad}\mathbf{P})\) by
Theorem~\ref{cor:dominio-realizacion-eliptica-compacta}. The inclusion
\[
 H^2(M,\operatorname{ad}\mathbf{P})
 \hookrightarrow
 L^2(M,\operatorname{ad}\mathbf{P})
\]
is compact by
Theorem~\ref{teo: rellich--kondrashov para haces vectoriales}. The
spectral theorem for operators with compact resolvent provides a constant
\(c_{\mathbf{A}}>0\) such that
\[
 \langle\Delta_{\mathbf{A}}^{(0)}\boldsymbol{\xi},\boldsymbol{\xi}\rangle_{L^2(M,\operatorname{ad}\mathbf{P})}
 \geq
 c_{\mathbf{A}}\|\boldsymbol{\xi}\|_{L^2(M,\operatorname{ad}\mathbf{P})}^2,
 \qquad
 \boldsymbol{\xi}\in H^2(M,\operatorname{ad}\mathbf P)\cap(\ker d_{\mathbf A})^{\perp_{L^2(M,\operatorname{ad}\mathbf P)}}.
\]

Let us make the construction of the inverse precise. Let
$(\boldsymbol\xi_j)$ be an orthonormal basis of eigensections, with
$\Delta_{\mathbf A}^{(0)}\boldsymbol\xi_j=\lambda_j\boldsymbol\xi_j$.
The kernel is $\ker d_{\mathbf A}$, since
$\langle\Delta_{\mathbf A}^{(0)}\boldsymbol\xi,\boldsymbol\xi\rangle_{L^2(M,\operatorname{ad}\mathbf P)}
=\|d_{\mathbf A}\boldsymbol\xi\|_{L^2(M,T^*M\otimes\operatorname{ad}\mathbf P)}^2$.
For $\mathbf f\in H^{s-1}(M,\operatorname{ad}\mathbf P)$ orthogonal
to this kernel, the series
\[
 \boldsymbol\xi
 =\sum_{\lambda_j>0}\lambda_j^{-1}
 \langle\mathbf f,\boldsymbol\xi_j\rangle_{L^2(M,\operatorname{ad}\mathbf P)}
 \boldsymbol\xi_j
\]
converges in $L^2(M,\operatorname{ad}\mathbf P)$, because
$\lambda_j\geq c_{\mathbf A}>0$ in the sum. Moreover,
$\displaystyle\sum_{j=1}^{\infty}\lambda_j^2|\langle\boldsymbol\xi,\boldsymbol\xi_j\rangle|^2
=\|\mathbf f\|_{L^2(M,\operatorname{ad}\mathbf P)}^2$;
the spectral characterization of the domain gives
$\Delta_{\mathbf A}^{(0)}\boldsymbol\xi=\mathbf f$ and
$\|\boldsymbol\xi\|_{L^2(M,\operatorname{ad}\mathbf P)}
\leq c_{\mathbf A}^{-1}\|\mathbf f\|_{L^2(M,\operatorname{ad}\mathbf P)}$.
A parametrix of order $-2$ and
Corollary~\ref{pseudo:cor-mapeo-sobolev-cerrada} give, for every real
$s\geq1$,
\[
 \|\boldsymbol\xi\|_{H^{s+1}(M,\operatorname{ad}\mathbf P)}
 \leq C_s\bigl(\|\mathbf f\|_{H^{s-1}(M,\operatorname{ad}\mathbf P)}
             +\|\boldsymbol\xi\|_{L^2(M,\operatorname{ad}\mathbf P)}\bigr)
 \leq C'_s\|\mathbf f\|_{H^{s-1}(M,\operatorname{ad}\mathbf P)}.
\]
The same energy identity proves injectivity on the complement
of the kernel. Thus
\[
 \Delta_{\mathbf{A}}^{(0)}\colon
 H^{s+1}(M,\operatorname{ad}\mathbf{P})\cap(\ker d_{\mathbf{A}})^\perp
 \longrightarrow
 H^{s-1}(M,\operatorname{ad}\mathbf{P})\cap(\ker d_{\mathbf{A}})^\perp
\]
is a topological isomorphism. If
\[
 \mathbf{a}\in H^s\left(
 M,T^*M\otimes\operatorname{ad}\mathbf{P}
 \right),
\]
then, for every \(\boldsymbol{\eta}\in\ker d_{\mathbf{A}}\),
\[
 \langle d_{\mathbf{A}}^*\mathbf{a},\boldsymbol{\eta}\rangle_{L^2(M,\operatorname{ad}\mathbf{P})}
 =
 \langle \mathbf{a},d_{\mathbf{A}}\boldsymbol{\eta}\rangle_{
 L^2(M,T^*M\otimes\operatorname{ad}\mathbf{P})}
 =
 0.
\]
Therefore \(d_{\mathbf{A}}^*\mathbf{a}\) belongs to the codomain of the preceding restriction.
There exists a unique
\[
 \boldsymbol{\xi}\in
 H^{s+1}(M,\operatorname{ad}\mathbf{P})\cap(\ker d_{\mathbf{A}})^\perp
\]
such that
\[
 \Delta_{\mathbf{A}}^{(0)}\boldsymbol{\xi}=d_{\mathbf{A}}^*\mathbf{a}.
\]
Set
\[
 \mathbf{b}:=\mathbf{a}-d_{\mathbf{A}}\boldsymbol{\xi}.
\]
Then \(d_{\mathbf{A}}^*\mathbf{b}=0\), so
\(\mathbf{a}=\mathbf{b}+d_{\mathbf{A}}\boldsymbol{\xi}\).

If \(d_{\mathbf{A}}\boldsymbol{\xi}\in\ker d_{\mathbf{A}}^*\), then
\[
 \|d_{\mathbf{A}}\boldsymbol{\xi}\|_{L^2(M,T^*M\otimes\operatorname{ad}\mathbf{P})}^2
 =
 \langle\Delta_{\mathbf{A}}^{(0)}\boldsymbol{\xi},\boldsymbol{\xi}\rangle_{L^2(M,\operatorname{ad}\mathbf{P})}
 =
 0,
\]
and hence \(d_{\mathbf{A}}\boldsymbol{\xi}=0\). The sum is direct. Finally, the
estimate for the inverse of \(\Delta_{\mathbf{A}}^{(0)}\) proves continuity of the
two projections.
\end{proof}

\begin{theorem}[Sobolev gauge group and compact Coulomb slice]
\label{prop:grupo-gauge-sobolev-compacto}
Let $(M,\mathbf{g})$ be a smooth compact connected Riemannian manifold of dimension $m$
without boundary; let $G$ be compact, let $\mathbf{P}\to M$ be a right principal
bundle, and let $A$ be a smooth connection. If $s>\frac{m}{2}$, then
\[
 \mathcal G^{s+1}(\mathbf{P})
 :=H^{s+1}(M,\operatorname{Ad}(\mathbf{P}))
\]
is a Hilbert Lie group modeled on
$H^{s+1}(M,\operatorname{ad}\mathbf{P})$, and its right action on
$\mathcal A^s(\mathbf{P})$ is smooth. The infinitesimal action is
\begin{equation}
\label{eq:accion-infinitesimal-gauge}
 \boldsymbol{\xi}\longmapsto d_A\boldsymbol{\xi}.
\end{equation}

Set
\[
 K_A:=\ker\bigl(d_A:H^{s+1}(M,\operatorname{ad}\mathbf{P})
 \longrightarrow H^s(M,T^*M\otimes\operatorname{ad}\mathbf{P})\bigr)
\]
and, using the $L^2(M,\operatorname{ad}\mathbf P)$ inner product induced by $\mathbf{g}$ and the
$\operatorname{Ad}$--invariant metric on $\operatorname{ad}\mathbf{P}$, define
\[\begin{gathered}
 X_A^{s+1}:=H^{s+1}(M,\operatorname{ad}\mathbf{P})\cap K_A^{\perp_{L^2(M,\operatorname{ad}\mathbf P)}},\\
 Y_A^{s-1}:=
 \left\{\mathbf{f}\in H^{s-1}(M,\operatorname{ad}\mathbf{P})\middle|
 \langle \mathbf{f},\boldsymbol{\zeta}\rangle_{\mathcal D',\mathcal D}=0
 \text{ for every }\boldsymbol{\zeta}\in K_A
 \right\}.
\end{gathered}\]
When $s\geq1$, the last space is precisely
$H^{s-1}(M,\operatorname{ad}\mathbf{P})\cap K_A^{\perp_{L^2(M,\operatorname{ad}\mathbf P)}}$; the annihilator formulation
remains meaningful if $s-1<0$, since $K_A$ consists of
smooth sections.
Then
\begin{equation}
\label{eq:laplaciano-gauge-isomorfismo-complementos}
 \Delta_A^{(0)}:=d_A^*d_A:
 X_A^{s+1}\xrightarrow{\ \cong\ }Y_A^{s-1}
\end{equation}
is a Hilbert space isomorphism and
\begin{equation}
\label{eq:descomposicion-coulomb-compacta}
 H^s(M,T^*M\otimes\operatorname{ad}\mathbf{P})
 =d_A X_A^{s+1}\oplus
 \ker\bigl(d_A^*:H^s(M,T^*M\otimes\operatorname{ad}\mathbf{P})
 \longrightarrow H^{s-1}(M,\operatorname{ad}\mathbf{P})\bigr).
\end{equation}

The stabilizer
\[
 \operatorname{Stab}(A)
 :=\{\mathbf{u}\in\mathcal G^{s+1}(\mathbf{P})\mid A\triangleleft \mathbf{u}=A\}
\]
is a compact finite-dimensional Lie group with Lie algebra $K_A$.
More precisely, if $p\in \mathbf{P}$ and $\operatorname{Hol}_p(A)\subseteq G$ is
the holonomy group based at $p$, evaluation at $p$ induces a
Lie group isomorphism
\begin{equation}
\label{eq:estabilizador-centralizador-holonomia}
 \operatorname{Stab}(A)
 \xrightarrow{\ \cong\ }
 Z_G\bigl(\operatorname{Hol}_p(A)\bigr).
\end{equation}

Finally, if
\[
 \mathcal S_A^s
 :=A+\ker\bigl(d_A^*:H^s(M,T^*M\otimes\operatorname{ad}\mathbf{P})
 \to H^{s-1}(M,\operatorname{ad}\mathbf{P})\bigr),
\]
there exist an open neighborhood $\mathcal U_A$ of $A$ in
$\mathcal A^s(\mathbf{P})$ and an open neighborhood $\mathcal W_A$ of $[A,e]$ in
$\mathcal S_A^s\times_{\operatorname{Stab}(A)}\mathcal G^{s+1}(\mathbf{P})$
for which the restriction of the map
\begin{equation}
\label{eq:tubo-local-rebanada-coulomb-compacta}
 \mathcal W_A\longrightarrow\mathcal U_A,
 \qquad [B,\mathbf{u}]\longmapsto B\triangleleft \mathbf{u},
\end{equation}
is a diffeomorphism. Here the equivalence relation is induced by
\[
 (B,\mathbf{u})\cdot \mathbf{h}=(B\triangleleft \mathbf{h},\mathbf{h}^{-1}\mathbf{u}),
 \qquad \mathbf{h}\in\operatorname{Stab}(A).
\]
Moreover, one can choose an open $\operatorname{Stab}(A)$-invariant neighborhood $\mathcal S_0$ of $A$ in
$\mathcal S_A^s$ such that
\[
 \mathcal S_0\times_{\operatorname{Stab}(A)}\mathcal G^{s+1}(\mathbf P)
 \longrightarrow\mathcal A^s(\mathbf P),
 \qquad [B,\mathbf u]\longmapsto B\triangleleft\mathbf u,
\]
is a diffeomorphism onto an invariant open set containing the orbit of
$A$. Every connection sufficiently close to $A$ has a representative in
$\mathcal S_0$, unique modulo $\operatorname{Stab}(A)$. Uniqueness refers
to representatives in this Coulomb neighborhood, rather than the entire
affine hyperplane $\mathcal S_A^s$.
\end{theorem}

\begin{proof}
Since $s>\frac{m}{2}$, we have $s+1>\frac{m}{2}+1$ and, in particular,
\[
H^{s+1}(M,\operatorname{ad}\mathbf{P})
\hookrightarrow C^0(M,\operatorname{ad}\mathbf{P};h_{\mathrm{ad}}).
\]
The Lie exponential is
$\operatorname{Ad}$--equivariant and therefore induces a bundle map
\[
 \operatorname{Exp}\colon\operatorname{ad}\mathbf{P}\longrightarrow
 \operatorname{Ad}(\mathbf{P}),
 \qquad [p,\xi]\longmapsto[p,\exp_G\xi].
\]
With respect to the $\operatorname{Ad}$--invariant inner product, choose an
$\operatorname{Ad}$--invariant ball $\mathcal O\subset\mathfrak g$ on
which $\exp_G$ is a diffeomorphism onto its image. For each section $\mathbf{u}$ of the
group bundle, the
formula
\[
 \Theta_{\mathbf{u}}(\boldsymbol{\xi}):=\mathbf{u}\,\operatorname{Exp}(\boldsymbol{\xi}),
 \qquad
 \|\boldsymbol{\xi}\|_{C^0(M,\operatorname{ad}\mathbf{P};h_{\mathrm{ad}})}<\operatorname{rad}(\mathcal O),
\]
has the pointwise inverse $v\mapsto\log_G(\mathbf{u}^{-1}v)$. In a finite cover
trivializing both bundles, the two formulas are superposition
operators. Lemma~\ref{lem:superposicion-sobolev-geometrica} proves that
they are smooth in $H^{s+1}(M,\operatorname{ad}\mathbf{P})$; equivalently, the differential at zero is the
identity, and Theorem~\ref{teo:funcion-inversa-local-banach} provides the
same local inverse. These charts, translated by each $\mathbf{u}$, construct the
Hilbert Lie group structure and the identification
$T_{\mathbf{u}}\mathcal G^{s+1}(\mathbf{P})
=\mathbf{u}H^{s+1}(M,\operatorname{ad}\mathbf{P})$.

In these charts, multiplication and inversion are superpositions of the
smooth finite-dimensional maps of $G$ and are therefore smooth. Likewise,
\[
\begin{gathered}
H^{s+1}(M,\operatorname{End}(\operatorname{ad}\mathbf{P}))
\cdot H^s(M,T^*M\otimes\operatorname{ad}\mathbf{P})
\longrightarrow H^s(M,T^*M\otimes\operatorname{ad}\mathbf{P}),\\[1ex]
\mathcal G^{s+1}(\mathbf{P})\longrightarrow H^s(M,T^*M\otimes\operatorname{ad}\mathbf{P}),
\qquad\mathbf{u}\longmapsto\mathbf{u}^{-1}d_A\mathbf{u}
\end{gathered}
\]
are smooth maps in charts. Here
$\mathbf u^{-1}d_A\mathbf u:=A\triangleleft\mathbf u-A$ is a global form
with values in $\operatorname{ad}\mathbf P$. The term
$\mathbf u^{-1}d\mathbf u$ alone is used in a local
trivialization: it combines with the connection potential to give that global form.
For real indices, apply
Corollary~\ref{cor:superposicion-hilbert-orden-real}. This constructs the Lie group and proves
smoothness of the action. Formula
\eqref{eq:diferencial-completo-accion-gauge-derecha}, with $\mathbf{u}=e$, gives
\eqref{eq:accion-infinitesimal-gauge}.

Proposition~\ref{prop:descomposicion-coulomb-compacta} and its proof
provide the spectral gap of $\Delta_A^{(0)}$, smoothness and finite
dimension of $K_A$, the isomorphism
\eqref{eq:laplaciano-gauge-isomorfismo-complementos}, and the decomposition
\eqref{eq:descomposicion-coulomb-compacta}, with continuous projections.
It applies directly when $s\geq1$.

For $0\leq s<1$, let $G_A$ be the reduced inverse of
$\Delta_A^{(0)}$, extended by zero on $K_A$. Its spectral construction
makes it self-adjoint on $L^2(M,\operatorname{ad}\mathbf P)$. The bound $G_A:L^2(M,\operatorname{ad}\mathbf P)\to H^2(M,\operatorname{ad}\mathbf P)$ gives
$G_A:H^{-2}(M,\operatorname{ad}\mathbf P)\to L^2(M,\operatorname{ad}\mathbf P)$ by duality; interpolation of these compatible
operators yields $G_A:H^{s-1}(M,\operatorname{ad}\mathbf P)\to H^{s+1}(M,\operatorname{ad}\mathbf P)$.
The projection $\Pi_{K_A}$ has finite-dimensional range consisting
of smooth sections and is therefore continuous at every Sobolev level.
Density of smooth sections gives
\[
 \Delta_A^{(0)}G_A=I-\Pi_{K_A}\text{ on }H^{s-1}(M,\operatorname{ad}\mathbf P),\qquad
 G_A\Delta_A^{(0)}=I-\Pi_{K_A}\text{ on }H^{s+1}(M,\operatorname{ad}\mathbf P).
\]
Moreover, $d_A^*\mathbf a$ annihilates $K_A$ by the duality identity
$\langle d_A^*\mathbf a,\boldsymbol\zeta\rangle
=\langle\mathbf a,d_A\boldsymbol\zeta\rangle=0$.
Thus $\boldsymbol\xi=G_A d_A^*\mathbf a$ and
$\mathbf b=\mathbf a-d_A\boldsymbol\xi$ give the same splitting and
continuous inverse in this range of indices. The section spaces
in these formulas are those of $\operatorname{ad}\mathbf P$ and
$T^*M\otimes\operatorname{ad}\mathbf P$ used in the statement.

We prove that every element of the Sobolev stabilizer is smooth. If $\mathbf{u}\in\mathcal G^{s+1}(\mathbf{P})$ satisfies
$A\triangleleft \mathbf{u}=A$, the local formula
\eqref{eq:accion-gauge-conexion} gives, in the weak sense,
\[
 \mathbf{u}^{-1}d\mathbf{u}=A-\operatorname{Ad}_{\mathbf{u}^{-1}}A.
\]
Since $s+1>\frac{m}{2}+1$, the section $\mathbf{u}$ is continuous. In a
trivialization of the group bundle over a coordinate domain $U$, with group coordinates in $\mathbb R^d$, $d=\dim G$, denote its local representation
by $\widehat u$. In group coordinates, the last equality is a system
$\partial_j\widehat u=F_j(x,\widehat u)$ with $F_j$ smooth. To gain one derivative, suppose
$\widehat u\in H^{s+1+j}_{\mathrm{loc}}(U,\mathbb R^d)$, with $j\in\mathbb N_0$.
Smooth superposition, applied to the map $(x,z)\mapsto F_i(x,z)$
in a relatively compact chart, gives
$F_i(x,\widehat u)\in H^{s+1+j}_{\mathrm{loc}}(U,\mathbb R^d)$ for each
$i\in\{1,\ldots,m\}$. The equation then provides all first
derivatives of $\widehat u$ in this space; the local Sobolev
characterization by derivatives implies
$\widehat u\in H^{s+2+j}_{\mathrm{loc}}(U,\mathbb R^d)$. Cutoff functions in slightly smaller
charts justify each step. Starting from $j=0$, this implication increases
regularity by one at each stage. A finite cover and Sobolev
embeddings give $\mathbf u\in C^\infty$. The equality is then classical, and
Theorem~\ref{teo:estabilizador-suave-centralizador-holonomia} shows that
$\mathbf{u}$ is a parallel section and evaluation induces the group
isomorphism \eqref{eq:estabilizador-centralizador-holonomia}. Conversely,
a smooth parallel section of $\operatorname{Ad}(\mathbf P)$ belongs to $\mathcal G^{s+1}(\mathbf P)$ because $M$ is compact,
so no element is lost.

Evaluation
$\operatorname{ev}_p:\mathcal G^{s+1}(\mathbf P)\to G$ is smooth:
in the chart $\Theta_{\mathbf u}$ it is the composition of continuous
linear evaluation of $\boldsymbol\xi\in H^{s+1}(M,\operatorname{ad}\mathbf P)$
with $\exp_G$ and multiplication by $\mathbf u(p)$. In the opposite direction, the extension of
$c\in Z_G(\operatorname{Hol}_p(A))$ is written in a trivialization as
$h(x)^{-1}c\,h(x)$, where $h(x)$ depends smoothly on the locally chosen parallel
transport. Thus the inverse is also smooth, and
\eqref{eq:estabilizador-centralizador-holonomia} is a Lie group
isomorphism.
Proposition~\ref{prop:centralizador-subgrupo-cerrado-lie} shows that the
centralizer is a closed Lie subgroup of the compact group $G$ and is therefore compact and finite-dimensional. Finally,
\eqref{eq:nucleo-dA-algebra-centralizador-holonomia} identifies its Lie algebra
with $K_A$.

Now apply the implicit function theorem. For small $\mathbf{a}$ and $\boldsymbol{\xi}$,
with $\boldsymbol{\xi}\in X_A^{s+1}$, set
\[
 \mathscr F(\mathbf{a},\boldsymbol{\xi})
 :=d_A^*\bigl((A+\mathbf{a})\triangleleft\exp_G(-\boldsymbol{\xi})-A\bigr)
 \in Y_A^{s-1}.
\]
The image belongs to $Y_A^{s-1}$ because $d_A^*$ takes values orthogonal
to $K_A$. The differential formula for the action gives
\begin{equation}
\label{eq:derivada-ift-coulomb-signo}
 D_{\boldsymbol{\xi}}\mathscr F(0,0)\boldsymbol{\eta}=-d_A^*d_A\boldsymbol{\eta}
 =-\Delta_A^{(0)}\boldsymbol{\eta},
\end{equation}
By
\eqref{eq:laplaciano-gauge-isomorfismo-complementos},
Theorem~\ref{teo:funcion-implicita-banach} produces a unique smooth function
$\mathbf{a}\mapsto\boldsymbol{\xi}(\mathbf{a})\in X_A^{s+1}$ such that
\[
 d_A^*\bigl((A+\mathbf{a})\triangleleft\exp_G(-\boldsymbol{\xi}(\mathbf{a}))-A\bigr)=0.
\]

Equivalently, the map
\[
 \Psi(\mathbf{b},\boldsymbol{\xi}):=(A+\mathbf{b})\triangleleft\exp_G(\boldsymbol{\xi})-A,
 \qquad
 d_A^*\mathbf{b}=0,\qquad \boldsymbol{\xi}\in X_A^{s+1},
\]
has differential
\[
 D\Psi_{(0,0)}(\mathbf{b},\boldsymbol{\xi})=\mathbf{b}+d_A\boldsymbol{\xi},
\]
which is an isomorphism by
\eqref{eq:descomposicion-coulomb-compacta}.
Theorem~\ref{teo:funcion-inversa-local-banach} gives the unique local parametrization
by $(\mathbf{b},\boldsymbol{\xi})$.

To incorporate the stabilizer, put
\[
 H:=\operatorname{Stab}(A),\qquad
 \mathfrak h:=K_A,\qquad X:=X_A^{s+1}.
\]
The affine Coulomb space $\mathcal S_A^s$ is $H$--invariant: if $\mathbf{h}$ is parallel,
then
\[
 (A+\mathbf{b})\triangleleft \mathbf{h}=A+\operatorname{Ad}_{\mathbf{h}^{-1}}\mathbf{b},
 \qquad
 d_A^*\operatorname{Ad}_{\mathbf{h}^{-1}}\mathbf{b}
 =\operatorname{Ad}_{\mathbf{h}^{-1}}d_A^*\mathbf{b}.
\]
The same identity shows that $\operatorname{Ad}_{\mathbf{h}}$ commutes with $d_A$.
Since $\mathbf{h}$ is parallel and the inner product on $\mathfrak g$ is
$\operatorname{Ad}$--invariant, $\operatorname{Ad}_{\mathbf{h}}$ commutes with the
covariant derivatives of the connection $A$. Choose the equivalent Sobolev
norm constructed with this connection; for real orders, we may use
powers of $I+\nabla_A^*\nabla_A$. The parallel endomorphism commutes
with this operator and its powers. It is therefore an
isometry of $H^{s+1}(M,\operatorname{ad}\mathbf{P})$. In particular, it preserves
both $K_A$ and its orthogonal complement $X$. Thus
\[
 H^{s+1}(M,\operatorname{ad}\mathbf{P})=\mathfrak h\oplus X
\]
is an $\operatorname{Ad}_H$--invariant decomposition.

By left translation,
the subbundle complementary to $TH$ in the restriction of the gauge
group tangent bundle to $H$ is
\[
 N_{\mathbf{h}}H:=(dL_{\mathbf{h}})_eX,
 \qquad
 NH\cong H\times X,
 \qquad (\mathbf{h},\boldsymbol{\xi})\longmapsto(dL_{\mathbf{h}})_e\boldsymbol{\xi}.
\]
Invariance of $X$ makes this identification compatible also
with right translations by elements of $H$. Consider
\[
 \Theta\colon H\times X\longrightarrow\mathcal G^{s+1}(\mathbf{P}),
 \qquad \Theta(\mathbf{h},\boldsymbol{\xi}):=\mathbf{h}\operatorname{Exp}(\boldsymbol{\xi}).
\]
This map is smooth by the Lie group structure already constructed. In the
left trivialization at $(\mathbf{h},0)$, its differential is
\[
 (dL_{\mathbf{h}^{-1}})_{\mathbf{h}}D\Theta_{(\mathbf{h},0)}
       \bigl((dL_{\mathbf{h}})_e\boldsymbol{\eta},\boldsymbol{\xi}\bigr)=\boldsymbol{\eta}+\boldsymbol{\xi},
 \qquad (\boldsymbol{\eta},\boldsymbol{\xi})\in\mathfrak h\oplus X.
\]
It is an isomorphism. The inverse function theorem gives, around each
$(\mathbf{h},0)$, an open product
$W_{\mathbf{h}}\times B_X(0,r_{\mathbf{h}})$ on which $\Theta$ is a diffeomorphism onto its
image, with $\mathbf{h}\in W_{\mathbf{h}}\subseteq H$ and $r_{\mathbf{h}}>0$. The open sets $W_{\mathbf{h}}$ cover
$H$. Choose a subcover
$W_{\mathbf{h}_1},\ldots,W_{\mathbf{h}_N}$ and set
\[
 \varepsilon_0:=\min_{1\leq i\leq N}r_{\mathbf{h}_i}>0.
\]
If $\mathbf{h}\in H$ and
$\|\boldsymbol{\xi}\|_{H^{s+1}(M,\operatorname{ad}\mathbf{P})}<\varepsilon_0$, some $W_{\mathbf{h}_i}$
contains $\mathbf{h}$ and
$(\mathbf{h},\boldsymbol{\xi})\in W_{\mathbf{h}_i}\times B_X(0,r_{\mathbf{h}_i})$. Thus $\Theta$ is a
local diffeomorphism on all of $H\times B_X(0,\varepsilon_0)$. Here and below,
$B_X$ denotes a ball for the norm
$H^{s+1}(M,\operatorname{ad}\mathbf{P})$ restricted to $X$.

Injectivity can be made uniform. If there were no
$0<\varepsilon\leq\varepsilon_0$ for which $\Theta$ was injective,
there would be sequences
\[
 \mathbf{h}_j,\mathbf{k}_j\in H,\qquad \boldsymbol{\xi}_j,\boldsymbol{\eta}_j\in X,\qquad
 \|\boldsymbol{\xi}_j\|_{H^{s+1}(M,\operatorname{ad}\mathbf{P})},
 \|\boldsymbol{\eta}_j\|_{H^{s+1}(M,\operatorname{ad}\mathbf{P})}
 \longrightarrow0,
\]
with $(\mathbf{h}_j,\boldsymbol{\xi}_j)\ne(\mathbf{k}_j,\boldsymbol{\eta}_j)$ and
$\Theta(\mathbf{h}_j,\boldsymbol{\xi}_j)=\Theta(\mathbf{k}_j,\boldsymbol{\eta}_j)$. After passing to subsequences,
compactness of $H$ gives $\mathbf{h}_j\to \mathbf{h}$ and $\mathbf{k}_j\to \mathbf{k}$. Passing to the limit
gives $\mathbf{h}=\mathbf{k}$. For large $j$, the two pairs then belong to a
common domain of a local inverse of $\Theta$ around $(\mathbf{h},0)$, forcing
$(\mathbf{h}_j,\boldsymbol{\xi}_j)=(\mathbf{k}_j,\boldsymbol{\eta}_j)$, a contradiction. Consequently,
for some $\varepsilon>0$,
\begin{equation}
\label{eq:tubo-estabilizador-gauge}
 \Theta\colon H\times B_X(0,\varepsilon)
 \xrightarrow{\ \cong\ }
 \mathcal N_\varepsilon:=
 \Theta\bigl(H\times B_X(0,\varepsilon)\bigr)
\end{equation}
is a diffeomorphism. Its image is open because $\Theta$ is a
local diffeomorphism and contains $H$. Thus
\eqref{eq:tubo-estabilizador-gauge} is a tubular neighborhood of
$H$, obtained from the complementary normal bundle $NH$. The terminology agrees with
Definition~\ref{def:variedades-riemannianas-variedad-4}: as in Fermi
coordinates, a neighborhood of the zero section of a normal bundle
is identified with a neighborhood of the submanifold; here the diffeomorphism is
$\Theta$ rather than a Riemannian exponential map on the Hilbert manifold.

The ball $B_X(0,\varepsilon)$ is $H$--invariant for the preceding norm. The
tube is equivariant under the right action
\[
 (\mathbf{h},\boldsymbol{\xi})\cdot \mathbf{k}=(\mathbf{h}\mathbf{k},\operatorname{Ad}_{\mathbf{k}^{-1}}\boldsymbol{\xi}),
 \qquad
 \Theta((\mathbf{h},\boldsymbol{\xi})\cdot \mathbf{k})=\Theta(\mathbf{h},\boldsymbol{\xi})\mathbf{k},
\]
and also under left translations. This is precisely the
compatibility required by the associated quotient.

To conclude, choose an open $H$--invariant neighborhood
\[
 \mathcal S_0\subset\mathcal S_A^s,\qquad A\in\mathcal S_0,
\]
and shrink $\mathcal S_0$ and $\varepsilon$ simultaneously so that
the local diffeomorphism $\Psi$ obtained above is injective on
$(\mathcal S_0-A)\times B_X(0,\varepsilon)$. Let $q$ be the projection onto the
associated quotient and set
\[
 \mathcal W_A:=q(\mathcal S_0\times\mathcal N_\varepsilon).
\]
We may take $\mathcal S_0$ to be a ball for an
$H$--invariant norm obtained by averaging any equivalent norm over the
compact group. The action of $H$ on
$\mathcal S_A^s\times\mathcal G^{s+1}(\mathbf{P})$ is free because of its second factor and
proper because $H$ is compact; in particular, the projection $q$ is
open and $\mathcal W_A$ is an open neighborhood of $[A,e]$.
If $\mathbf{u}=\mathbf{h}\operatorname{Exp}(\boldsymbol{\xi})$ belongs to $\mathcal N_\varepsilon$,
then, in the associated quotient,
\[
 [B,\mathbf{h}\operatorname{Exp}(\boldsymbol{\xi})]
 =[B\triangleleft \mathbf{h},\operatorname{Exp}(\boldsymbol{\xi})],
\]
with $B\triangleleft \mathbf{h}\in\mathcal S_0$ after shrinking
$\mathcal S_0$ invariantly. This representative is unique.
Indeed, if
\[
 [B_1,\operatorname{Exp}(\boldsymbol{\xi}_1)]
 =[B_2,\operatorname{Exp}(\boldsymbol{\xi}_2)],
\]
there exists $\mathbf{k}\in H$ such that
$B_2=B_1\triangleleft \mathbf{k}$ and
$\operatorname{Exp}(\boldsymbol{\xi}_2)=\mathbf{k}^{-1}\operatorname{Exp}(\boldsymbol{\xi}_1)$. Hence
\[
 \Theta(e,\boldsymbol{\xi}_1)=\Theta(\mathbf{k},\boldsymbol{\xi}_2),
\]
and injectivity of the tube implies $\mathbf{k}=e$, $\boldsymbol{\xi}_1=\boldsymbol{\xi}_2$, and $B_1=B_2$.
Thus
\[
 \mathcal S_0\times B_X(0,\varepsilon)
 \longrightarrow\mathcal W_A,
 \qquad (B,\boldsymbol{\xi})\longmapsto[B,\operatorname{Exp}(\boldsymbol{\xi})]
\]
is a smooth bijective chart with smooth inverse induced by
\eqref{eq:tubo-estabilizador-gauge}. In this chart, the map
\eqref{eq:tubo-local-rebanada-coulomb-compacta} agrees with
\[
 (B,\boldsymbol{\xi})\longmapsto B\triangleleft\operatorname{Exp}(\boldsymbol{\xi}),
\]
that is, the parametrization $\Psi$ translated by $A$. Its image
$\mathcal U_A$ is open, and the restriction is a diffeomorphism. This proves the local parametrization. To obtain a slice for
the full group, we must also control gauge transformations
that have not been chosen in advance in $\mathcal N_\varepsilon$.

Fix a faithful unitary representation of $G$ with associated vector bundle
$\mathbf E$. Gauge transformations are represented as
unitary endomorphisms of $\mathbf E$, and $d_A\mathbf u$ here denotes their
covariant derivative before multiplication by $\mathbf u^{-1}$.
The elliptic estimate for the fixed connection $A$ gives
\[
 \|\mathbf u\|_{H^{s+1}(M,\operatorname{End}\mathbf E)}
 \leq C_A\bigl(\|d_A\mathbf u\|_{H^s(M,T^*M\otimes\operatorname{End}\mathbf E)}
                  +\|\mathbf u\|_{L^2(M,\operatorname{End}\mathbf E)}\bigr).
\]
This estimate applies to arbitrary endomorphisms, not just group
elements. If $(A+\mathbf b_1)\triangleleft\mathbf u=A+\mathbf b_2$,
the action formula is equivalent to
\[
 d_A\mathbf u=\mathbf u\mathbf b_2-\mathbf b_1\mathbf u.
\]
Sobolev product estimates imply
\[
 \|d_A\mathbf u\|_{H^s(M,T^*M\otimes\operatorname{End}\mathbf E)}
 \leq C_{\mathrm{mult}}\bigl(\|\mathbf b_1\|_{H^s(M,T^*M\otimes\operatorname{ad}\mathbf P)}+\|\mathbf b_2\|_{H^s(M,T^*M\otimes\operatorname{ad}\mathbf P)}\bigr)
          \|\mathbf u\|_{H^{s+1}(M,\operatorname{End}\mathbf E)}.
\]
Take fixed $C_A,C_{\mathrm{mult}}>0$; they depend on $A$, $s$, the
metric, and the representation. The constant $C_{\mathrm{mult}}$ includes
the norm of the induced homomorphism
$\operatorname{ad}\mathbf P\to\operatorname{End}\mathbf E$. If $d$ is the complex rank of $\mathbf E$,
the Hilbert--Schmidt norm of a unitary endomorphism is $\sqrt d$, hence
\[
 \|\mathbf u\|_{L^2(M,\operatorname{End}\mathbf E)}
       =\sqrt{d\,\lambda_{\mathbf g}(M)}.
\]
Choose $\delta>0$ with $\delta\leq(4C_A C_{\mathrm{mult}})^{-1}$.
For $\|\mathbf b_i\|_{H^s(M,T^*M\otimes\operatorname{ad}\mathbf P)}<\delta$,
$i\in\{1,2\}$, the preceding two estimates imply
\[
 C_A C_{\mathrm{mult}}
 \bigl(\|\mathbf b_1\|_{H^s(M,T^*M\otimes\operatorname{ad}\mathbf P)}+\|\mathbf b_2\|_{H^s(M,T^*M\otimes\operatorname{ad}\mathbf P)}\bigr)\leq\frac12,
 \qquad
 \|\mathbf u\|_{H^{s+1}(M,\operatorname{End}\mathbf E)}
       \leq 2C_A\sqrt{d\,\lambda_{\mathbf g}(M)}.
\]
Indeed, the term containing the $H^{s+1}(M,\operatorname{End}\mathbf E)$ norm of $\mathbf u$ on the
right has coefficient at most $1/2$ and can be moved to the left.
The resulting bound is independent of $\mathbf u$ and
$\mathbf b_1,\mathbf b_2$ in these balls.

Suppose there exist $\mathbf b_{1,j},\mathbf b_{2,j}\to0$ in $H^s(M,T^*M\otimes\operatorname{ad}\mathbf P)$
and transformations $\mathbf u_j$ relating the two connections but
remaining outside $\mathcal N_\varepsilon$. The preceding bound and Rellich's theorem
allow extraction of a subsequence converging in $H^s(M,\operatorname{End}\mathbf E)$ and, by embedding, in
$C^0$, to a group-valued section $\mathbf u$. The differential identity
shows that $d_A\mathbf u_j\to0$ in $H^s(M,T^*M\otimes\operatorname{End}\mathbf E)$, so
$d_A\mathbf u=0$ in distributions. The stabilizer regularity argument
proves that $\mathbf u\in H=\operatorname{Stab}(A)$.
Applying the elliptic estimate to $\mathbf u_j-\mathbf u$ gives
\[
 \|\mathbf u_j-\mathbf u\|_{H^{s+1}(M,\operatorname{End}\mathbf E)}
 \leq C_A\bigl(\|d_A\mathbf u_j\|_{H^s(M,T^*M\otimes\operatorname{End}\mathbf E)}
                      +\|\mathbf u_j-\mathbf u\|_{L^2(M,\operatorname{End}\mathbf E)}\bigr)\longrightarrow0.
\]
This contradicts the fact that $\mathcal N_\varepsilon$ is an open neighborhood of
$H$. Consequently, after shrinking $\mathcal S_0$, every transformation
relating two points of $\mathcal S_0$ belongs to
$\mathcal N_\varepsilon$.

Write this transformation as
$\mathbf u=\mathbf h\operatorname{Exp}\boldsymbol\xi$, with
$\mathbf h\in H$ and $\boldsymbol\xi\in X$. If
$B_1\triangleleft\mathbf u=B_2$, with $B_1,B_2\in\mathcal S_0$,
invariance under $H$ keeps $B_1\triangleleft\mathbf h$ in
$\mathcal S_0$. Injectivity of the parametrization $\Psi$ compares the
pairs $(B_1\triangleleft\mathbf h,\boldsymbol\xi)$ and $(B_2,0)$ and gives
$\boldsymbol\xi=0$, $B_2=B_1\triangleleft\mathbf h$. This proves uniqueness
with respect to any group element.

Finally, the differential of the parametrization remains invertible
at $(B,e)$ for every $B$ in a sufficiently small neighborhood of $A$,
by continuity and openness of the set of isomorphisms. Translation by
the action transfers this property to $(B,\mathbf u)$ for any
$\mathbf u$. The map from the associated quotient is thus a local diffeomorphism
at all its points, and the preceding argument proves its injectivity.
Its image is open and invariant and contains the orbit of $A$; being an
open injective map with smooth local inverses, it is a diffeomorphism
onto that image. This completes the slice construction.
\end{proof}

The gain of one derivative in the gauge element is visible both in
$\mathbf{u}^{-1}d\mathbf{u}$ and in the elliptic isomorphism used in the Coulomb argument.

On a compact base, $d_A^*d_A$ has compact resolvent: its spectrum is
discrete, each eigenspace is finite-dimensional, and there is no essential
spectrum. The elliptic estimate on $(\ker d_A)^{\perp}$ replaces, in this
case, the additional gap assumption arising on a noncompact base.

On a complete noncompact base without boundary, we say that a smooth
connection $B$ satisfies the connection condition $(B_k)$ if there exist constants
$D_0,\ldots,D_k>0$ such that
\begin{equation}
\label{eq:geometria-acotada-conexion-gauge}
 \|(\nabla^B)^i\mathbf F_B\|_\infty\leq D_i,
 \qquad 0\leq i\leq k.
\end{equation}
For fixed $A_0$, define its affine Sobolev component by
\[
\mathcal A^r_{A_0}(\mathbf{P})
:=
A_0+H^r(M,T^*M\otimes\operatorname{ad}\mathbf{P}).
\]
On a noncompact base, a gauge transformation is likewise not required to
be integrable as a function with values in $G$. Fix a faithful unitary
representation of $G$ and use the connection induced by $A_0$. Consider
sections $\mathbf u$ of $\operatorname{Ad}(\mathbf P)$ of class
$H^{r+1}_{\mathrm{loc}}(M,\operatorname{Ad}\mathbf P)$ for which
\[
 \mathbf u^{-1}d_{A_0}\mathbf u
 :=A_0\triangleleft\mathbf u-A_0
 \in H^r(M,T^*M\otimes\operatorname{ad}\mathbf P).
\]
This set is $\mathcal G^{r+1}(A_0)$. Its topology is described by the
charts $\boldsymbol\xi\mapsto\mathbf u\operatorname{Exp}\boldsymbol\xi$,
with $\boldsymbol\xi\in H^{r+1}(M,\operatorname{ad}\mathbf P)$ of small
uniform norm. The norms are always based on $A_0$ and $\mathbf g$.

The definition is compatible with multiplication. If
$\mathbf b_{\mathbf u}=\mathbf u^{-1}d_{A_0}\mathbf u$, the action law gives
\[
 \mathbf b_{\mathbf u\mathbf v}
   =\operatorname{Ad}_{\mathbf v^{-1}}\mathbf b_{\mathbf u}
      +\mathbf b_{\mathbf v},
 \qquad
 \mathbf b_{\mathbf u^{-1}}=-\operatorname{Ad}_{\mathbf u}\mathbf b_{\mathbf u}.
\]
The coefficients $\mathbf u$ are bounded because $G$ is compact.
Estimates for their derivatives follow from $d_{A_0}\mathbf u=\mathbf u\mathbf b_{\mathbf u}$:
in the step from order $j$ to $j+1$, Leibniz's rule
differentiates one of the factors $\mathbf b_{\mathbf u}$ or adds a new factor
by differentiating $\mathbf u$. For $j\in\{0,\ldots,r\}$, Sobolev
products, with $r>m/2$, control all terms. In particular,
$\operatorname{Ad}_{\mathbf u}$ and its inverse are multipliers in $H^r(M,T^*M\otimes\operatorname{ad}\mathbf P)$.
The preceding two identities preserve the membership condition, and
the charts construct the Hilbert Lie group described in
\cite[\S3]{EichhornHeberGauge}. Its action on the affine component is written as
\[
 (A_0+\mathbf a)\triangleleft\mathbf u
 =A_0+\mathbf b_{\mathbf u}+\operatorname{Ad}_{\mathbf u^{-1}}\mathbf a.
\]

Let $\mathcal G_e^{r+1}(A_0)$ be the identity component of this group.
A central constant section $c\in Z(G)$ has $d_{A_0}c=0$, so it
belongs to the full group even if it does not belong to that component.
If $Z(G)$ is countable, set
\[
 \mathcal G_\bullet^{r+1}(A_0)
   :=\bigcup_{c\in Z(G)}c\,\mathcal G_e^{r+1}(A_0).
\]
If the center is uncountable, set
$\mathcal G_\bullet^{r+1}(A_0)=\mathcal G_e^{r+1}(A_0)$.
The selected components form a subgroup. Since the Hilbert model
is separable, the identity component is second countable; a
countable union of its translates preserves this property. This is the
restriction used in the slice theorem.

Essential spectrum and a gap above the kernel are understood as in
Definition~\ref{def:espectro-esencial-brecha-cap24}.
Proposition~\ref{prop:brecha-espectral-rango-inversa} proves the
equivalence of gap, closed range, and bounded inverse, and justifies
the passage from a positive essential spectral bound for the restriction.

\begin{proposition}[Regularity and Coulomb decomposition under a spectral gap]
\label{prop:coulomb-sobolev-brecha-uniforme}
Let $(M,\mathbf g)$ be a smooth Riemannian manifold of dimension $m$ without
boundary, of bounded geometry of order $k$. Let $\mathbf P\to M$ be a
principal bundle with compact structure group $G$, and let $A_0$ be a smooth connection
satisfying \eqref{eq:geometria-acotada-conexion-gauge}. Fix integers
$k-1\geq r>m/2+2$ and a connection
$A=A_0+\mathbf a$, where
$\mathbf a\in H^r(M,T^*M\otimes\operatorname{ad}\mathbf P)$.
In this proposition we write
\[
 \begin{gathered}
 E:=\operatorname{ad}\mathbf P,\qquad F:=T^*M\otimes E,\\
 H^j(M,E):=H^j(M,E;\nabla^{A_0}),\qquad H^j(M,F):=H^j(M,F;\nabla^{A_0}).
 \end{gathered}
\]
Norms and adjoints use $\mathbf g$ and an
$\operatorname{Ad}(G)$--invariant fiber metric. Let $\Delta_A^{(0)}$
be the nonnegative operator associated with the closed form
\[
 q_A(\boldsymbol\xi,\boldsymbol\eta)
 :=\int_M\langle d_A\boldsymbol\xi,d_A\boldsymbol\eta\rangle
                \,d\lambda_{\mathbf g},
 \qquad \boldsymbol\xi,\boldsymbol\eta\in H^1(M,E),
\]
and set $K_A:=\ker\Delta_A^{(0)}\subseteq L^2(M,E)$.
Then $\mathcal D(\Delta_A^{(0)})=H^2(M,E)$. For each $R>0$ there is
$C_{r,R}>0$, depending only on $r,R$ and the fixed geometric data,
such that, if $\|\mathbf a\|_{H^r(M,F)}\leq R$,
\begin{equation}
\label{eq:regularidad-gauge-coeficientes-sobolev}
 \begin{gathered}
 \boldsymbol\xi\in\mathcal D(\Delta_A^{(0)}),\quad
 \Delta_A^{(0)}\boldsymbol\xi\in H^{r-1}(M,E)
 \quad\Longrightarrow\quad \boldsymbol\xi\in H^{r+1}(M,E),\\
 \|\boldsymbol\xi\|_{H^{r+1}(M,E)}
 \leq C_{r,R}\bigl(\|\Delta_A^{(0)}\boldsymbol\xi\|_{H^{r-1}(M,E)}
                       +\|\boldsymbol\xi\|_{L^2(M,E)}\bigr).
 \end{gathered}
\end{equation}
If, in addition, there is a gap $\gamma>0$ above $K_A$, the inverse
$G_A$ of $\Delta_A^{(0)}$ on $K_A^{\perp_{L^2(M,\operatorname{ad}\mathbf P)}}$ satisfies
\begin{equation}
\label{eq:inversa-gauge-brecha-sobolev}
 \|G_A\mathbf f\|_{H^{r+1}(M,E)}
 \leq C_{r,R}(1+\gamma^{-1})\|\mathbf f\|_{H^{r-1}(M,E)},
 \qquad \mathbf f\in H^{r-1}(M,E)\cap K_A^{\perp_{L^2(M,\operatorname{ad}\mathbf P)}}.
\end{equation}
Consequently,
\begin{equation}
\label{eq:coulomb-sobolev-brecha-suma}
 \begin{gathered}
 H^r(M,F)=Y_A\oplus d_A X_A,\\
 Y_A:=\ker(d_A^*:H^r(M,F)\to H^{r-1}(M,E)),\qquad
 X_A:=H^{r+1}(M,E)\cap K_A^{\perp_{L^2(M,\operatorname{ad}\mathbf P)}}.
 \end{gathered}
\end{equation}
The sum is orthogonal in $L^2(M,F)$ and its projections are continuous
on $H^r(M,F)$.
\end{proposition}

\begin{proof}
Write $d_0$ for the connection induced by $A_0$ on $E$, and set
$T_{\mathbf a}\boldsymbol\xi=[\mathbf a,\boldsymbol\xi]$.
Lemma~\ref{lem:modulo-sobolev-uniforme-entero}, also applied to
the contractions defining the adjoint, gives constants $c_j>0$ such that
\[
 \begin{aligned}
 \|T_{\mathbf a}\boldsymbol\xi\|_{H^j(M,F)}
 &\leq c_j R\|\boldsymbol\xi\|_{H^j(M,E)},\\
 \|T_{\mathbf a}^*\mathbf v\|_{H^j(M,E)}
 &\leq c_j R\|\mathbf v\|_{H^j(M,F)},\qquad 0\leq j\leq r.
 \end{aligned}
\]
In particular, $d_A=d_0+T_{\mathbf a}$ in distributions.
The norm $\bigl(\|\boldsymbol\xi\|_{L^2(M,E)}^2+
q_A(\boldsymbol\xi,\boldsymbol\xi)\bigr)^{1/2}$ is equivalent to the
$H^1(M,E)$ norm: the two inequalities follow from
$d_A=d_0+T_{\mathbf a}$ and $d_0=d_A-T_{\mathbf a}$, using
$\|T_{\mathbf a}\|_{L^2(M,E)\to L^2(M,F)}\leq c_0R$.
Thus the form is closed and densely defined.
Theorem~\ref{teo:representacion-formas-cerradas} constructs its
self-adjoint operator. By density of smooth compactly supported sections,
its weak equation is $d_A^*d_A\boldsymbol\xi=\mathbf f$.

We specify the estimate for the smooth reference operator
$\Delta_0=d_0^*d_0$ used in the regularity argument. For each integer
$1\leq j\leq r$,
\begin{equation}
\label{eq:regularidad-referencia-orden-finito-gauge}
 \|\mathbf v\|_{H^{j+1}(M,E)}
 \leq C_j\bigl(\|\Delta_0\mathbf v\|_{H^{j-1}(M,E)}
                    +\|\mathbf v\|_{H^j(M,E)}\bigr),
\end{equation}
provided $\mathbf v\in H^j(M,E)$ and $\Delta_0\mathbf v\in H^{j-1}(M,E)$
in distributions. Only the geometric bounds through the finite orders
fixed here are needed. To verify this, use the normal charts,
synchronous frames and cutoffs of
Proposition~\ref{prop:covariantes-localizacion-orden-entero}.
In these charts the operator is
\[
 L_i=-\sum_{\alpha,\beta=1}^m g_i^{\alpha\beta}
                         \partial_\alpha\partial_\beta
       +\sum_{\alpha=1}^m B_i^\alpha\partial_\alpha+C_i.
\]
Bounded geometry of order $k\geq r+1$ provides a positive ellipticity
constant and uniform bounds for the coefficients and their derivatives
through order $r$. The coefficient freezing argument leading to
\eqref{eq:estimacion-interior-B-W2p-cap13}, with $p=2$, uses only
ellipticity, a bound for the first derivatives of
$g_i^{\alpha\beta}$ and bounds for $B_i^\alpha,C_i$.
It therefore applies to $L_i$ and gives the second order estimate on
two concentric balls of fixed radii, with a constant independent of $i$.

To obtain order $j+1$, apply this estimate to
$\partial^\nu\mathbf v$ for $|\nu|\leq j-1$, using intermediate
balls between the preceding two. Each term in
$[L_i,\partial^\nu]\mathbf v$ contains a coefficient derivative
of order between $1$ and $|\nu|$, multiplied by a derivative of
$\mathbf v$ of order at most $|\nu|+1\leq j$.
Cutoff terms likewise contain derivatives of
$\mathbf v$ of order at most $j$. There are finitely many terms and
intermediate radii, depending only on $j$; their bounds give
\[
 \|\mathbf v\|_{H^{j+1}(V_i,E)}
 \leq C_j'\bigl(\|\Delta_0\mathbf v\|_{H^{j-1}(U_i,E)}
                       +\|\mathbf v\|_{H^j(U_i,E)}\bigr).
\]
Local elliptic regularity in
Corollary~\ref{pseudo:cor-regularidad-eliptica-sobolev-local} allows
this inequality to be applied to distributional solutions, complexifying
the bundle if it is real: it first
gives $H^{j+1}_{\mathrm{loc}}(M,\operatorname{ad}\mathbf P)$, and approximation on each larger ball
extends the estimates obtained for smooth sections.
Squaring and summing over $i\in I$, bounded multiplicity of
$(U_i)_{i\in I}$ and equivalence of localized norms prove
\eqref{eq:regularidad-referencia-orden-finito-gauge}. This argument
also proves finiteness of the global norm on the left.

For $\boldsymbol\xi\in\mathcal D(\Delta_A^{(0)})$, write
$\mathbf f=\Delta_A^{(0)}\boldsymbol\xi$. In distributions,
\begin{equation}
\label{eq:expansion-laplaciano-conexion-sobolev}
 \Delta_0\boldsymbol\xi
 =\mathbf f-d_0^*(T_{\mathbf a}\boldsymbol\xi)
             -T_{\mathbf a}^*d_0\boldsymbol\xi
             -T_{\mathbf a}^*T_{\mathbf a}\boldsymbol\xi.
\end{equation}
If $\boldsymbol\xi\in H^j(M,E)$, $1\leq j\leq r$, the three additional
terms belong to $H^{j-1}(M,E)$ and their sum has norm at most
$C_j''(R+R^2)\|\boldsymbol\xi\|_{H^j(M,E)}$.
We start with $\boldsymbol\xi\in H^1(M,E)$, the form domain.
If $\mathbf f\in L^2(M,E)$, the preceding identity and
\eqref{eq:regularidad-referencia-orden-finito-gauge}, with $j=1$,
give $\boldsymbol\xi\in H^2(M,E)$. Conversely, for a section in
$H^2(M,E)$, the expression $d_A^*d_A\boldsymbol\xi$ belongs to $L^2(M,\operatorname{ad}\mathbf P)$;
integration by parts and density in $H^1(M,E)$ identify it with
the form operator. Its domain is therefore exactly $H^2(M,E)$.

If $\mathbf f\in H^{r-1}(M,E)$, repeat this step for
$j=1,\ldots,r$. The constants satisfy the recurrence
\[
 \|\boldsymbol\xi\|_{H^{j+1}(M,E)}
 \leq C_j\|\mathbf f\|_{H^{j-1}(M,E)}
   +C_j\bigl(1+C_j''(R+R^2)\bigr)
                    \|\boldsymbol\xi\|_{H^j(M,E)}.
\]
The initial norm is controlled by the form identity:
\[
 \begin{aligned}
 \|d_A\boldsymbol\xi\|_{L^2(M,F)}^2
 &=\langle\mathbf f,\boldsymbol\xi\rangle_{L^2(M,E)},\\
 \|\boldsymbol\xi\|_{H^1(M,E)}
 &\leq \tfrac12\|\mathbf f\|_{L^2(M,E)}
             +(\tfrac32+c_0R)\|\boldsymbol\xi\|_{L^2(M,E)}.
 \end{aligned}
\]
The second inequality uses
$\sqrt{xy}\leq(x+y)/2$, for $x,y\geq0$.
The recurrence has $r$ steps and proves
\eqref{eq:regularidad-gauge-coeficientes-sobolev}, with a constant
uniform for $\|\mathbf a\|_{H^r(M,F)}\leq R$.
In particular, $K_A\subseteq H^{r+1}(M,E)$ and $K_A=\ker d_A$ in $L^2(M,\operatorname{ad}\mathbf P)$,
where the latter equality is understood on the form domain.

Now let $\mathbf f\in H^{r-1}(M,E)\cap K_A^{\perp_{L^2(M,\operatorname{ad}\mathbf P)}}$.
The gap and the spectral theorem give
$\boldsymbol\xi=G_A\mathbf f\in\mathcal D(\Delta_A^{(0)})
\cap K_A^{\perp_{L^2(M,\operatorname{ad}\mathbf P)}}$, with
$\|\boldsymbol\xi\|_{L^2(M,E)}\leq\gamma^{-1}
\|\mathbf f\|_{L^2(M,E)}$. Substitution into
\eqref{eq:regularidad-gauge-coeficientes-sobolev} yields
\eqref{eq:inversa-gauge-brecha-sobolev}.

The algebraic construction of
Proposition~\ref{prop:descomposicion-coulomb-compacta} now applies
with this inverse. Indeed, if $\mathbf v\in H^r(M,F)$ and
$\boldsymbol\eta\in K_A$, density in $H^1(M,E)$ gives
\[
 \langle d_A^*\mathbf v,\boldsymbol\eta\rangle_{L^2(M,E)}
 =\langle\mathbf v,d_A\boldsymbol\eta\rangle_{L^2(M,F)}=0.
\]
Thus the following operators are defined:
\[
 Q_A\mathbf v:=G_A d_A^*\mathbf v\in X_A,\qquad
 \Pi_A\mathbf v:=\mathbf v-d_AQ_A\mathbf v\in Y_A.
\]
If $D_A,D_A^*>0$ bound, respectively, the norms of
$d_A:H^{r+1}(M,E)\to H^r(M,F)$ and
$d_A^*:H^r(M,F)\to H^{r-1}(M,E)$, we may take
\begin{equation}
\label{eq:cotas-proyecciones-coulomb-brecha}
 \|Q_A\|\leq C_{r,R}(1+\gamma^{-1})D_A^*,\qquad
 \|\Pi_A\|\leq1+D_A C_{r,R}(1+\gamma^{-1})D_A^*.
\end{equation}
The identities $d_A^*\Pi_A=0$, $Q_A d_A\boldsymbol\xi=\boldsymbol\xi$
for $\boldsymbol\xi\in X_A$, and $\Pi_A+d_AQ_A=I$ show that
$\Pi_A$ and $d_AQ_A$ are the projections of
\eqref{eq:coulomb-sobolev-brecha-suma}.
The image of $d_A$ on all of $H^{r+1}(M,E)$ equals $d_AX_A$:
the orthogonal projection onto $K_A$ maps $L^2(M,E)$ into $H^{r+1}(M,E)$
by \eqref{eq:regularidad-gauge-coeficientes-sobolev}, applied with
$\mathbf f=0$, and subtracting it does not change the image.
Finally, the integration by parts identity above proves orthogonality.
Both summands are closed as images of continuous projections.
\end{proof}

\pagebreak
\begin{proposition}[Persistence of the essential spectrum and the gap]
\label{prop:persistencia-espectral-componente-gauge}
Under the hypotheses and notation of
Proposition~\ref{prop:coulomb-sobolev-brecha-uniforme}, every
$A=A_0+\mathbf a\in\mathcal A^r_{A_0}(\mathbf P)$ satisfies
\begin{equation}
\label{eq:invariancia-espectro-esencial-gauge}
 \sigma_{\mathrm e}(\Delta_A^{(0)})
 =\sigma_{\mathrm e}(\Delta_{A_0}^{(0)}).
\end{equation}
If $M$ is also connected and noncompact, then
$\ker\Delta_A^{(0)}=\{0\}$ for every $A$ in this component. In this case,
the condition
$\displaystyle\inf\sigma_{\mathrm e}(\Delta_{A_0}^{(0)})>0$
implies that each $A$ admits a constant $\gamma_A>0$ such that
\[
 \|d_A\boldsymbol\xi\|_{L^2(M,F)}^2
 \geq\gamma_A\|\boldsymbol\xi\|_{L^2(M,E)}^2,
 \qquad \boldsymbol\xi\in H^1(M,E).
\]
There is a constant $c\geq1$, depending only on the bound for the
fiberwise bracket and on the embedding
$H^r(M,F)\hookrightarrow C_b^0(M,F)$, with the following property: if
$B\in\mathcal A^r_{A_0}(\mathbf P)$ and
\[
 \|B-A\|_{H^r(M,F)}<\frac{\sqrt{\gamma_A}}{2c},
\]
then
\begin{equation}
\label{eq:brecha-local-uniforme-gauge}
 \|d_B\boldsymbol\xi\|_{L^2(M,F)}^2
 \geq\frac{\gamma_A}{4}\|\boldsymbol\xi\|_{L^2(M,E)}^2,
 \qquad \boldsymbol\xi\in H^1(M,E).
\end{equation}
\end{proposition}

\begin{proof}
Write $\Delta_0=\Delta_{A_0}^{(0)}$ and
$T_{\mathbf a}\boldsymbol\xi=[\mathbf a,\boldsymbol\xi]$.
Identity~\eqref{eq:expansion-laplaciano-conexion-sobolev} gives
\[
 K_{\mathbf a}:=\Delta_A^{(0)}-\Delta_0
 =d_0^*T_{\mathbf a}+T_{\mathbf a}^*d_0
      +T_{\mathbf a}^*T_{\mathbf a}
 \colon H^2(M,E)\longrightarrow L^2(M,E).
\]
In the first term, the Leibniz rule gives contractions of
$\nabla^{A_0}\mathbf a$ with $\boldsymbol\xi$ and of $\mathbf a$ with
$\nabla^{A_0}\boldsymbol\xi$. The second contains only products of
the latter kind, and the third contains two factors $\mathbf a$ and
one factor $\boldsymbol\xi$. Hence there is $C>0$, depending only on
the bounds for these contractions and the uniform embedding into
$C_b^1(M,T^*M\otimes\operatorname{ad}\mathbf P)$, such that
\begin{equation}
\label{eq:continuidad-perturbacion-laplaciano-gauge}
 \|(K_{\mathbf a}-K_{\mathbf b})\boldsymbol\xi\|_{L^2(M,E)}
 \leq C\bigl(1+\|\mathbf a\|_{H^r(M,F)}+\|\mathbf b\|_{H^r(M,F)}\bigr)
       \|\mathbf a-\mathbf b\|_{H^r(M,F)}\|\boldsymbol\xi\|_{H^2(M,E)}.
\end{equation}
Here $\mathbf a,\mathbf b\in H^r(M,F)$ and $\boldsymbol\xi\in H^2(M,E)$.
For the quadratic term we use
\[
 T_{\mathbf a}^*T_{\mathbf a}-T_{\mathbf b}^*T_{\mathbf b}
 =T_{\mathbf a-\mathbf b}^*T_{\mathbf a}
       +T_{\mathbf b}^*T_{\mathbf a-\mathbf b}.
\]

Choose $\mathbf a_n\in\Gamma_c(F)$, $n\in\mathbb N$, such that
$\mathbf a_n\to\mathbf a$ in $H^r(M,F)$. This density follows from
the cutoffs in \eqref{eq:cola-densidad-Wkp-geometria-acotada-cap13}
and the Meyers--Serrin Theorem~\ref{meyers-serrin-haz}.
For each $n$, choose a relatively compact domain $D_n\subset M$
with smooth boundary containing $\operatorname{supp}\mathbf a_n$, and
$\chi_n\in C_c^\infty(D_n)$ equal to one on a neighborhood of that support.
Then $K_{\mathbf a_n}\boldsymbol\xi
=K_{\mathbf a_n}(\chi_n\boldsymbol\xi)$. This operator factors as
\[
 H^2(M,E)\xrightarrow{\ \boldsymbol\xi\mapsto\chi_n\boldsymbol\xi\ }
 H^2(D_n,E)\hookrightarrow H^1(D_n,E)
 \xrightarrow{\ K_{\mathbf a_n}\ }L^2(M,E),
\]
where the last map includes extension by zero of the result, whose
support is contained in $\operatorname{supp}\mathbf a_n$.
The first and last arrows are continuous, and the middle inclusion
is compact by Corollary~\ref{cor:rellich-descenso-un-orden-haces}.
Thus $K_{\mathbf a_n}:H^2(M,E)\to L^2(M,E)$ is compact.
Estimate~\eqref{eq:continuidad-perturbacion-laplaciano-gauge} gives
convergence in operator norm to $K_{\mathbf a}$, so $K_{\mathbf a}$
is compact as well.

The domains of $\Delta_A^{(0)}$ and $\Delta_0$ are both $H^2(M,E)$.
The step $j=1$ in the proof of
Proposition~\ref{prop:coulomb-sobolev-brecha-uniforme} gives, for some
$C_A>0$,
\[
 \|\boldsymbol\xi\|_{H^2(M,E)}
 \leq C_A\bigl(\|\Delta_A^{(0)}\boldsymbol\xi\|_{L^2(M,E)}
                     +\|\boldsymbol\xi\|_{L^2(M,E)}\bigr),
 \qquad\boldsymbol\xi\in H^2(M,E).
\]
The reverse estimate follows from the continuity of
$\Delta_A^{(0)}:H^2(M,E)\to L^2(M,E)$. Its graph norm is therefore
equivalent to the $H^2(M,\operatorname{ad}\mathbf P)$ norm; the same holds for $\Delta_0$.
For each $\lambda\in\mathbb R$, the operators
$\Delta_A^{(0)}-\lambda I$ and $\Delta_0-\lambda I$, viewed as maps from
$H^2(M,E)$ to $L^2(M,E)$, differ by the compact operator $K_{\mathbf a}$.
Theorem~\ref{teo:estabilidad-indice-fredholm}, applied in both directions,
shows that one is Fredholm if and only if the other is. Definition
\ref{def:espectro-esencial-brecha-cap24} yields
\eqref{eq:invariancia-espectro-esencial-gauge}.

Now assume that $M$ is connected and noncompact. Choose
$\rho_0>0$ less than one fourth of the injectivity radius.
Bounded geometry provides $v_0>0$ such that
$\lambda_{\mathbf g}(B_{\mathbf g}(x,\rho_0))\geq v_0$ for every $x\in M$.
Completeness and noncompactness provide infinitely many disjoint balls
of this radius: after finitely many centers have been chosen, the union
of the closed balls of radius $2\rho_0$ about them is compact by
Hopf--Rinow and does not cover $M$, so the next center can be chosen
outside that union. Hence $\lambda_{\mathbf g}(M)=\infty$.
If $\boldsymbol\xi\in\ker\Delta_A^{(0)}$, the form identity gives
$d_A\boldsymbol\xi=0$, and the regularity already established gives
$\boldsymbol\xi\in H^{r+1}(M,E)\subset C^1(M,E)$.
Metric compatibility implies that $|\boldsymbol\xi|$ is constant on
the connected base. This constant is zero because
$\boldsymbol\xi\in L^2(M,E)$ and the volume is infinite.

Equality of the essential spectra and
Proposition~\ref{prop:brecha-espectral-rango-inversa} now provide
$\gamma_A>0$. The spectral theorem extends the operator inequality
to its form domain $H^1(M,E)$, since
$q_A(\boldsymbol\xi,\boldsymbol\xi)
=\|(\Delta_A^{(0)})^{1/2}\boldsymbol\xi\|_{L^2(M,E)}^2$.
Finally, let $c_{\mathrm{ad}}$ be a pointwise bound for
$|[\mathbf v,\boldsymbol\xi]|/(|\mathbf v||\boldsymbol\xi|)$, and let
$C_{\mathrm{emb}}$ be the constant of the embedding
$H^r(M,F)\hookrightarrow C_b^0(M,F)$.
Set $c=\max\{1,c_{\mathrm{ad}}C_{\mathrm{emb}}\}$. For every
$\boldsymbol\xi\in H^1(M,E)$,
\[
 \begin{aligned}
 \|d_B\boldsymbol\xi\|_{L^2(M,F)}
 &\geq\|d_A\boldsymbol\xi\|_{L^2(M,F)}
                -c\|B-A\|_{H^r(M,F)}\|\boldsymbol\xi\|_{L^2(M,E)}\\
 &\geq\bigl(\sqrt{\gamma_A}-c\|B-A\|_{H^r(M,F)}\bigr)
                                      \|\boldsymbol\xi\|_{L^2(M,E)}.
 \end{aligned}
\]
On the stated ball the last factor in parentheses is at least
$\sqrt{\gamma_A}/2$. Squaring yields
\eqref{eq:brecha-local-uniforme-gauge}.
\end{proof}

\pagebreak
\begin{corollary}[Local parametrization in Coulomb gauge]
\label{cor:carta-coulomb-sobolev-brecha}
Under the hypotheses and notation of
Proposition~\ref{prop:coulomb-sobolev-brecha-uniforme}, assume that
the gap $\gamma>0$ exists. There are $\rho,\varepsilon>0$ such that
\[
 \Psi_A(\mathbf b,\boldsymbol\xi)
 :=(A+\mathbf b)\triangleleft\operatorname{Exp}\boldsymbol\xi
\]
is a diffeomorphism from
$B_{Y_A\times X_A}(0,\rho)$ onto an open neighborhood of $A$ containing
$A+B_{H^r(M,F)}(0,\varepsilon)$. On the product we use the norm
$\|\mathbf b\|_{H^r(M,F)}+\|\boldsymbol\xi\|_{H^{r+1}(M,E)}$.
Thus each connection in the latter ball has a unique representation
of this form with $(\mathbf b,\boldsymbol\xi)$ in the stated domain.
\end{corollary}

\begin{proof}
The derivative at $(0,0)$ is the isomorphism
\[
 L_A(\mathbf b,\boldsymbol\xi)=\mathbf b+d_A\boldsymbol\xi,
 \qquad L_A^{-1}\mathbf v=(\Pi_A\mathbf v,Q_A\mathbf v).
\]
Fix $B_A>0$ such that $\|L_A^{-1}\|\leq B_A$, using
\eqref{eq:cotas-proyecciones-coulomb-brecha}, and $c\geq1$ such that
\[
 \|[\boldsymbol\xi,\mathbf v]\|_{H^r(M,F)}
 \leq c\|\boldsymbol\xi\|_{H^{r+1}(M,E)}\|\mathbf v\|_{H^r(M,F)}.
\]
The action formula and differentiation of the exponential give
\begin{equation}
\label{eq:carta-coulomb-serie-controlada}
 \Psi_A(\mathbf b,\boldsymbol\xi)-A
 =e^{-\operatorname{ad}\boldsymbol\xi}\mathbf b
   +\int_0^1e^{-t\operatorname{ad}\boldsymbol\xi}
                               d_A\boldsymbol\xi\,dt.
\end{equation}
The integral is a Bochner integral in $H^r(M,F)$.
The identity is first obtained for smooth sections by differentiating
$\exp(-t\boldsymbol\xi)d_A\exp(t\boldsymbol\xi)$ with respect to $t$.
The module estimates and approximation in $H^{r+1}(M,E)$ and
$H^r(M,F)$ extend it to the stated spaces. The exponential series converge
there, together with all their derivatives, on bounded balls:
the term of order $\ell$ has operator norm at most
$(c\|\boldsymbol\xi\|_{H^{r+1}(M,E)})^\ell/\ell!$
on $H^r(M,F)$. The same bound, with the integer factors arising from
differentiation, ensures convergence of the differentiated series.
In particular, the map is smooth.

Set $F_A:=\Psi_A-A$ and $K:=c(2+2D_A)e^c$.
For $0<\rho\leq1$ and $\|z\|_{Y_A\times X_A}<\rho$, the preceding
formula implies
\begin{equation}
\label{eq:cota-diferencial-carta-coulomb}
 \|DF_A(z)-L_A\|\leq K\rho.
\end{equation}
To verify the constant, the derivative with respect to $\mathbf b$
differs from the identity by at most $c\rho e^c$ in norm.
The derivative of the first term with respect to $\boldsymbol\xi$
contributes at most $c\rho e^c$; after subtracting
$d_A\boldsymbol\eta$, the derivative of the second contributes at most
$cD_A\rho e^c+cD_A\rho e^c/2$ per unit of
$\|\boldsymbol\eta\|_{H^{r+1}(M,E)}$.
Here we use $e^t-1\leq te^t$, for $t\geq0$, and
$\displaystyle\int_0^1t\,dt=1/2$. The chosen value of $K$ dominates
both bounds for the sum norm on the product.

Let $C_{\mathrm{emb}}$ be an embedding constant for
$H^{r+1}(M,E)\hookrightarrow C^0_b(M,E)$, and let $\iota_G>0$ be a radius
on which the exponential of $G$ is injective. This radius can be chosen
invariant under the adjoint action because $G$ is compact. Choose
\[
 0<\rho\leq
 \min\{1,\iota_G/(2C_{\mathrm{emb}}),(2B_AK)^{-1}\},
 \qquad \varepsilon:=\rho/(4B_A).
\]
For $R_A(z):=L_A^{-1}F_A(z)-z$, estimate
\eqref{eq:cota-diferencial-carta-coulomb} gives
$\|DR_A(z)\|\leq1/2$ on the ball of radius $\rho$, and $R_A(0)=0$.
The fundamental theorem of calculus along segments implies
\[
 \|R_A(z)-R_A(w)\|\leq\tfrac12\|z-w\|,
 \qquad z,w\in B_{Y_A\times X_A}(0,\rho).
\]
Hence $F_A$ is injective on this ball. Its derivative is also invertible
at each point by the Neumann series applied to $I+DR_A(z)$;
the inverse function theorem shows that the image is open and the
inverse is smooth.
For $\|\mathbf v\|_{H^r(M,F)}<\varepsilon$, the map
\[
 z\longmapsto L_A^{-1}\mathbf v-R_A(z)
\]
maps the closed ball of radius $\rho/2$ into itself and is a contraction
with constant at most $1/2$. Its unique fixed point satisfies
$F_A(z)=\mathbf v$ and $\|z\|\leq2B_A\|\mathbf v\|_{H^r(M,F)}$.
This proves the asserted inclusion and the unique representation
within the stated domain.
\end{proof}

\begin{lemma}[Control of transformations in the identity component]
\label{lem:control-global-calibre-brecha}
Under the assumptions of
Proposition~\ref{prop:coulomb-sobolev-brecha-uniforme}, suppose in addition
that $M$ is connected and noncompact and that the spectral gap holds at $A$.
Retain $E=\operatorname{ad}\mathbf P$, $F=T^*M\otimes E$, and the
norms based on $A_0$. There exist $J_A\geq1$ and $\delta_A>0$ such that,
for $\mathbf u\in\mathcal G_e^{r+1}(A_0)$, the condition
\[
 \boldsymbol\beta:=A\triangleleft\mathbf u-A,\qquad
 \|\boldsymbol\beta\|_{H^r(M,F)}<\delta_A
\]
implies $\mathbf u=\operatorname{Exp}\boldsymbol\xi$, with
$\boldsymbol\xi\in H^{r+1}(M,E)$ in the domain of the pointwise logarithm and
\begin{equation}
\label{eq:control-logaritmo-global-calibre}
 \|\boldsymbol\xi\|_{H^{r+1}(M,E)}
 \leq J_A\|\boldsymbol\beta\|_{H^r(M,F)}.
\end{equation}
In particular, the stabilizer of $A$ in
$\mathcal G_e^{r+1}(A_0)$ is trivial.
\end{lemma}

\begin{proof}
Proposition~\ref{prop:persistencia-espectral-componente-gauge}
gives infinite volume and $K_A=\{0\}$. Consequently,
$Q_Ad_A=I$ on $H^{r+1}(M,E)$.
With $C_A:=\max\{1,\|Q_A\|\}$ we obtain
\begin{equation}
\label{eq:coercividad-primer-orden-calibre}
 \|\boldsymbol\zeta\|_{H^{r+1}(M,E)}
 \leq C_A\|d_A\boldsymbol\zeta\|_{H^r(M,F)},
 \qquad\boldsymbol\zeta\in H^{r+1}(M,E),
\end{equation}
with the bound for $Q_A$ in
\eqref{eq:cotas-proyecciones-coulomb-brecha}.

In the faithful unitary representation fixed above,
$\mathbf u-I\in H^{r+1}(M,\operatorname{End}\mathbf E)$ for every
$\mathbf u\in\mathcal G_e^{r+1}(A_0)$.
Indeed, an exponential neighborhood of the identity generates an open
and closed subgroup containing the identity component. Thus each
$\mathbf u$ is a finite product of exponentials of
sections in $H^{r+1}(M,E)$. The identity
$\displaystyle e^{\boldsymbol\zeta}-I=
\int_0^1e^{t\boldsymbol\zeta}\boldsymbol\zeta\,dt$ and the unitized
Sobolev algebra prove the assertion for each factor and
for their product.

Fix a normal radius $\iota_G>0$ for a bi-invariant metric on $G$
and $0<\eta<\iota_G/4$. Let $\chi\colon[0,\infty)\to[0,1]$ be smooth,
equal to $1$ on $[0,4\eta^2]$ and to $0$ on $[9\eta^2,\infty)$.
The map
\[
 f(g):=\chi(|\log g|^2)\log g
\]
on the normal ball of radius $\iota_G$, extended by zero outside it,
is smooth and equivariant under conjugation, and satisfies $f(e)=0$.
Hence $f(\mathbf u)$ is a section of $E$.
Compactness of $G$ bounds the derivatives of $f$ through order $r+1$.

The action identity is $d_A\mathbf u=\mathbf u\boldsymbol\beta$.
The first derivative of $f(\mathbf u)$ introduces a factor
$\boldsymbol\beta$. Each subsequent derivative differentiates an existing
factor or a coefficient depending on $\mathbf u$; in the latter case it adds
another factor $\boldsymbol\beta$. Equivariance allows this
rule to be applied with the connections induced by $A$.
The identities for weak derivatives follow by local
approximation in $H^{r+1}(M,\operatorname{Ad}\mathbf P)$ and the product estimates.
Induction gives, for $j\in\{1,\ldots,r+1\}$,
\[
 |(\nabla^A)^j f(\mathbf u)|
 \leq D_j\sum_{q=1}^{j}
 \sum_{\substack{\alpha_1,\ldots,\alpha_q\geq0\\
                 \alpha_1+\cdots+\alpha_q=j-q}}
 \prod_{\nu=1}^{q}|(\nabla^A)^{\alpha_\nu}\boldsymbol\beta|
 \quad\text{almost everywhere on }M.
\]
The indices $\alpha_\nu$ are integers; the products represent the
norms of the tensor products in the Leibniz rules.
The constants $D_j>0$ bound the coefficients in those rules and the
derivatives of $f$ on $G$, and do not depend on $\mathbf u$.

The covariant norms for $A$ and $A_0$ through order $r+1$ are
equivalent. Expanding
$\nabla^A=\nabla^{A_0}+[A-A_0,\cdot]$, the Sobolev module estimate controls
the lower order terms; the same identity, solved for
$\nabla^{A_0}$, gives the reverse estimate. The constants depend on
$r$, the geometric bounds, and $\|A-A_0\|_{H^r(M,F)}$.
In each product above, the total order is $j-q\leq r$.
The Sobolev multiplication theorem bounds its $L^2(M,(T^*M)^{\otimes j}\otimes E)$ norm by a
constant times $\|\boldsymbol\beta\|_{H^r(M,F)}^q$.
Summing the orders $j=1,\ldots,r+1$ gives a constant
$C_{A,f}\geq1$ such that
\begin{equation}
\label{eq:logaritmo-truncado-calibre}
 \|d_A f(\mathbf u)\|_{H^r(M,F)}
 \leq C_{A,f}\sum_{q=1}^{r+1}\|\boldsymbol\beta\|_{H^r(M,F)}^q.
\end{equation}
Moreover, $|f(g)|\leq C_f|g-I|$ on $G$ for some constant $C_f$.
This also controls the term of order zero and gives
$f(\mathbf u)\in H^{r+1}(M,E)$ before applying
\eqref{eq:coercividad-primer-orden-calibre}.

Let $C_{\mathrm{emb}}>0$ be a constant for the embedding
$H^{r+1}(M,E)\hookrightarrow C_b^0(M,E)$. Set
\[
 J_A:=(r+1)C_AC_{A,f},\qquad
 \delta_A:=\min\{1/2,\eta/(2C_{\mathrm{emb}}J_A)\}.
\]
The preceding estimates imply
\[
 \|f(\mathbf u)\|_{H^{r+1}(M,E)}
 \leq J_A\|\boldsymbol\beta\|_{H^r(M,F)},\qquad
 \|f(\mathbf u)\|_\infty<\eta/2.
\]
There is a point $x_0\in M$ where $\mathbf u(x_0)$ is at distance less than
$\eta/2$ from the identity. Otherwise, compactness of $G$ would give
a uniform positive lower bound for $|\mathbf u-I|$, contradicting
its integrability. Distance to the identity is well defined in
the adjoint bundle because the metric is invariant under conjugation.
If the continuous section $\mathbf u$ left the normal ball of radius
$\eta$, a path from $x_0$ to a point outside would have
a first point with $|\log\mathbf u|=\eta$. There
$|f(\mathbf u)|=\eta$, contradicting the preceding bound.
Thus $\boldsymbol\xi=f(\mathbf u)=\log\mathbf u$ throughout $M$,
and \eqref{eq:control-logaritmo-global-calibre} holds.
If $A\triangleleft\mathbf u=A$, that estimate gives
$\boldsymbol\xi=0$ and $\mathbf u=I$.
\end{proof}

\begin{proposition}[Gauge tube over the orbit]
\label{prop:tubo-calibre-brecha-conexa}
Under the assumptions of Lemma~\ref{lem:control-global-calibre-brecha},
the orbit $\mathcal O_A:=A\triangleleft\mathcal G_\bullet^{r+1}(A_0)$
is a closed submanifold of $\mathcal A^r_{A_0}(\mathbf P)$ and admits
an equivariant tube constructed from $Y_A=\ker d_A^*\subset H^r(M,F)$.
Write $S(A)$ for the stabilizer of $A$ in the full group
$\mathcal G^{r+1}(A_0)$.
There exists $\varepsilon_A>0$ such that
\[
 \mathcal S_A:=\{A+\mathbf b\mid\mathbf b\in Y_A,
                               \|\mathbf b\|_{r,A}<\varepsilon_A\}
\]
is a slice invariant under the full stabilizer $S(A)$ and
\begin{equation}
\label{eq:producto-global-rebanada-calibre}
 \begin{aligned}
 \Theta_A\colon\mathcal S_A\times\mathcal G_e^{r+1}(A_0)
      &\longrightarrow\mathcal A^r_{A_0}(\mathbf P),\\
 (B,\mathbf u)&\longmapsto B\triangleleft\mathbf u
 \end{aligned}
\end{equation}
is a diffeomorphism onto an open invariant neighborhood of
$\mathcal O_A$. The norm $\|\cdot\|_{r,A}$ is defined in the proof.
\end{proposition}

\begin{proof}
For $B\in\mathcal A^r_{A_0}(\mathbf P)$ and $\mathbf v\in H^r(M,F)$,
define
\begin{equation}
\label{eq:metrica-fuerte-invariante-calibre}
 \|\mathbf v\|_{r,B}^2
 :=\sum_{j=0}^{r}\int_M|(\nabla^B)^j\mathbf v|^2\,d\lambda_{\mathbf g}.
\end{equation}
Its polarization is a smooth strong Riemannian metric. The
finite expansions of the covariant derivatives and the Sobolev module
estimate show that its coefficients are continuous polynomial maps of
$B-A_0$; the preceding norm equivalence proves that it is strong.
The identity
\[
 (\nabla^{B\triangleleft\mathbf u})^j
 (\operatorname{Ad}_{\mathbf u^{-1}}\mathbf v)
 =\operatorname{Ad}_{\mathbf u^{-1}}(\nabla^B)^j\mathbf v,
 \qquad j\in\{0,\ldots,r\},
\]
is obtained first for smooth data and then by local approximation
and the product estimates. Together with invariance of the inner product
on the adjoint bundle, it proves gauge invariance of the metric.

Fix $\Lambda_A\geq1$ such that
\begin{equation}
\label{eq:comparacion-metrica-fuerte-calibre}
 \Lambda_A^{-1}\|\mathbf v\|_{H^r(M,F)}
 \leq\|\mathbf v\|_{r,B}\leq\Lambda_A\|\mathbf v\|_{H^r(M,F)},
 \qquad\|B-A\|_{H^r(M,F)}<1.
\end{equation}
The constant can be chosen uniformly on this ball, since
$\|B-A_0\|_{H^r(M,F)}\leq\|A-A_0\|_{H^r(M,F)}+1$.
Let $d_r$ be the infimum of the lengths of piecewise
$C^1$ paths for this metric. The segment from $A$ to
$A+\mathbf b$ gives
$d_r(A,A+\mathbf b)\leq\Lambda_A\|\mathbf b\|_{H^r(M,F)}$ if
$\|\mathbf b\|_{H^r(M,F)}<1$.
If $0<s<1$, every path reaching the outside of the ball of
radius $s$ has length at least $s/\Lambda_A$: apply the
lower bound up to its first exit. Consequently,
\begin{equation}
\label{eq:bola-distancia-calibre}
 d_r(A,B)<s/\Lambda_A
 \quad\Longrightarrow\quad\|B-A\|_{H^r(M,F)}<s.
\end{equation}
These bounds around each connection prove that $d_r$ is a distance
inducing the Sobolev topology. The gauge action acts by isometries.

The countable center of a compact Lie group is finite. Put
$C_\bullet=Z(G)$ if the center is countable and $C_\bullet=\{e\}$ otherwise.
Then
$\mathcal G_\bullet^{r+1}(A_0)=C_\bullet\mathcal G_e^{r+1}(A_0)$.
The first factor acts trivially. The intersection of the factors is
trivial, since a constant central section different from $I$ does not
differ from $I$ by an $L^2(M,\operatorname{End}\mathbf E)$ section. The lemma gives
\begin{equation}
\label{eq:estabilizador-grupo-seleccionado-calibre}
 S(A)\cap\mathcal G_\bullet^{r+1}(A_0)=C_\bullet.
\end{equation}

Let $\rho>0$ be the radius from
Corollary~\ref{cor:carta-coulomb-sobolev-brecha}. Set
\[
 s_A:=\min\{1/2,\delta_A/2,\rho/(4J_A)\},\qquad
 \varepsilon_A:=\min\{\rho/(4\Lambda_A),1/(4\Lambda_A),
                                      s_A/(4\Lambda_A^3)\}.
\]
Suppose $B_i=A+\mathbf b_i\in\mathcal S_A$ for $i=1,2$ and
$B_1\triangleleft\mathbf u=B_2$, with
$\mathbf u\in\mathcal G_\bullet^{r+1}(A_0)$. By invariance,
\[
 d_r(A,A\triangleleft\mathbf u)
 \leq d_r(A,B_2)+d_r(B_1,A)
 <2\Lambda_A^2\varepsilon_A<s_A/\Lambda_A.
\]
Thus $\|A\triangleleft\mathbf u-A\|_{H^r(M,F)}<s_A$.
Writing $\mathbf u=c\mathbf v$ with $c\in C_\bullet$ and
$\mathbf v\in\mathcal G_e^{r+1}(A_0)$, the lemma gives
$\mathbf v=\operatorname{Exp}\boldsymbol\xi$ and
$\|\boldsymbol\xi\|_{H^{r+1}(M,E)}<\rho/4$.
Also $\|\mathbf b_i\|_{H^r(M,F)}<\rho/4$. The equality
$\Psi_A(\mathbf b_1,\boldsymbol\xi)=\Psi_A(\mathbf b_2,0)$ lies
within the injectivity ball of the corollary. Hence
$\boldsymbol\xi=0$, $\mathbf b_1=\mathbf b_2$, and $\mathbf u=c$.
This controls every element of the selected group.

The differential of $\Theta_A$ at $(B,I)$ is that of $\Psi_A$ at
$(B-A,0)$ and is invertible. The action and group translations
transport this property to $(B,\mathbf u)$ for every $\mathbf u$.
Thus $\Theta_A$ is a local diffeomorphism. If
$B_1\triangleleft\mathbf u_1=B_2\triangleleft\mathbf u_2$, the preceding
argument applied to $\mathbf u_1\mathbf u_2^{-1}$ proves
injectivity. The image $\mathcal U_A$ is open and invariant, and the
local inverses agree, so $\Theta_A$ is a diffeomorphism
onto it. For $\mathbf h\in S(A)$, the action on $A+\mathbf b$
is $A+\operatorname{Ad}_{\mathbf h^{-1}}\mathbf b$ and preserves
$\ker d_A^*$ and $\|\cdot\|_{r,A}$ by the covariant identities.
Thus $\mathcal S_A$ is invariant under the full stabilizer,
which may be larger than $C_\bullet$.

The orbit corresponds to $\{A\}\times\mathcal G_e^{r+1}(A_0)$ and is a
closed submanifold of $\mathcal U_A$. It is closed in the whole component:
if $B\in\overline{\mathcal O_A}\setminus\mathcal O_A$, choose
$A\triangleleft\mathbf u_n\to B$ and define
$B_n=B\triangleleft\mathbf u_n^{-1}$. Then
$B_n\in\overline{\mathcal O_A}\setminus\mathcal O_A$ and
\[
 d_r(A,B_n)=d_r(A\triangleleft\mathbf u_n,B)\longrightarrow0.
\]
Thus $B_n\to A$ in the Sobolev topology and, for all sufficiently
large $n$, $B_n\in\mathcal U_A$. This contradicts closedness of the orbit
in that neighborhood.

The normal complement with respect to the $L^2(M,T^*M\otimes\operatorname{ad}\mathbf P)$ inner product is the smooth bundle
\[
 \mathcal N_A:=\{(C,\mathbf w)\mid C\in\mathcal O_A,
                   \mathbf w\in H^r(M,F),\ d_C^*\mathbf w=0\}.
\]
Its trivialization is
$(\mathbf b,\mathbf u)\mapsto
(A\triangleleft\mathbf u,\operatorname{Ad}_{\mathbf u^{-1}}\mathbf b)$,
with $\mathbf b\in Y_A$ and $\mathbf u\in\mathcal G_e^{r+1}(A_0)$.
Invariance shows that
\[
 E_A(C,\mathbf w):=C+\mathbf w,\qquad
 (C,\mathbf w)\in\mathcal N_A,\quad\|\mathbf w\|_{r,C}<\varepsilon_A,
\]
is the preceding diffeomorphism expressed in this bundle. The radial compression
\[
 (C,\mathbf w)\longmapsto
 \left(C,\frac{\varepsilon_A\mathbf w}
                  {\sqrt{1+\|\mathbf w\|_{r,C}^2}}\right)
\]
is a diffeomorphism from the whole bundle onto the indicated normal disk;
its inverse on each fiber is
\mbox{$\mathbf z\mapsto\mathbf z/\sqrt{\varepsilon_A^2-\|\mathbf z\|_{r,C}^2}$}.
Composition with $E_A$ gives the equivariant tube.
To place it inside a prescribed open invariant neighborhood of the
orbit, reduce $\varepsilon_A$ until the slice lies in that neighborhood;
invariance then ensures that the whole tube lies there.
\end{proof}

\begin{corollary}[Structure of the quotient by the selected components]
\label{cor:cociente-hilbert-brecha-gauge}
Under the hypotheses of
Proposition~\ref{prop:persistencia-espectral-componente-gauge}, assume
that $M$ is connected and noncompact and that
$\displaystyle\inf\sigma_{\mathrm e}(\Delta_{A_0}^{(0)})>0$.
Then
\[
 \mathcal R:=\mathcal A^r_{A_0}(\mathbf P)/\mathcal G_\bullet^{r+1}(A_0)
 =\mathcal A^r_{A_0}(\mathbf P)/\mathcal G_e^{r+1}(A_0)
\]
is a connected, Hausdorff, second countable Hilbert manifold.
The projection $\pi:\mathcal A^r_{A_0}(\mathbf P)\to\mathcal R$ is a
smooth submersion. If $d_r$ is the length distance associated with
\eqref{eq:metrica-fuerte-invariante-calibre}, the formula
\begin{equation}
\label{eq:distancia-cociente-calibre}
 \overline d_r([A],[B])
 :=\inf_{\mathbf u\in\mathcal G_e^{r+1}(A_0)}
                      d_r(A,B\triangleleft\mathbf u)
\end{equation}
defines a distance inducing the quotient topology.
\end{corollary}

\begin{proof}
Constant central transformations act trivially on connections.
The description
$\mathcal G_\bullet^{r+1}(A_0)=C_\bullet\mathcal G_e^{r+1}(A_0)$
in Proposition~\ref{prop:tubo-calibre-brecha-conexa} proves that the
two quotients agree. Spectral persistence allows us to apply that
proposition at every point of the affine component. In particular,
all orbits are closed, and we have the diffeomorphisms
\[
 \Theta_A:\mathcal S_A\times\mathcal G_e^{r+1}(A_0)
                  \longrightarrow\mathcal U_A
\]
of \eqref{eq:producto-global-rebanada-calibre}.

The distance $d_r$ is finite: the affine segment between two
connections has finite length because the strong metric varies
continuously along the compact parameter interval. Moreover, every
gauge transformation is an isometry. These facts show that
\eqref{eq:distancia-cociente-calibre} is independent of the
representatives and is symmetric. To check the triangle inequality,
let $A,B,C$ be three connections and let
$\mathbf u,\mathbf v\in\mathcal G_e^{r+1}(A_0)$.
The right action gives
\[
 \begin{aligned}
 \overline d_r([A],[C])
 &\leq d_r(A,C\triangleleft(\mathbf v\mathbf u))\\
 &\leq d_r(A,B\triangleleft\mathbf u)
            +d_r(B\triangleleft\mathbf u,
                         (C\triangleleft\mathbf v)\triangleleft\mathbf u)\\
 &=d_r(A,B\triangleleft\mathbf u)+d_r(B,C\triangleleft\mathbf v).
 \end{aligned}
\]
Take the infimum first over $\mathbf u$ and then over $\mathbf v$.
If $\overline d_r([A],[B])=0$, there are
$\mathbf u_n\in\mathcal G_e^{r+1}(A_0)$, $n\in\mathbb N$, such that
$d_r(A,B\triangleleft\mathbf u_n)<1/n$.
Since $d_r$ induces the Sobolev topology and the orbit of $B$ is closed,
$A$ belongs to that orbit. Thus $[A]=[B]$, so $\overline d_r$ is a distance.

The projection $\pi$ is open, since the preimage of the image of an
open set is the union of its translates by the group. For every $A$
and every $\varepsilon>0$,
\[
 B_{\overline d_r}([A],\varepsilon)
                =\pi\bigl(B_{d_r}(A,\varepsilon)\bigr).
\]
This identity and the invariance of saturated open sets show that
the topology of $\overline d_r$ is the quotient topology. More
explicitly, if $\pi^{-1}(V)$ is open and $[A]\in V$, some ball
$B_{d_r}(A,\varepsilon)$ is contained in $\pi^{-1}(V)$, and its image
is the ball above, contained in $V$.
The affine component is connected and has a countable basis of open
sets, since it is modeled on the separable Hilbert space $H^r(M,F)$.
Its image is connected, and the images of this basis form a countable
basis for the quotient.

Each restriction $\pi|_{\mathcal S_A}$ is a homeomorphism onto the
open set $\pi(\mathcal U_A)$. Indeed, $\Theta_A$ identifies the orbits
in $\mathcal U_A$ with the fibers of the projection onto $\mathcal S_A$.
We therefore obtain charts
\[
 \chi_A:\pi(\mathcal U_A)\longrightarrow
             \{\mathbf b\in Y_A\mid\|\mathbf b\|_{r,A}<\varepsilon_A\},
 \qquad \chi_A([A+\mathbf b])=\mathbf b.
\]
Let us check their smooth compatibility. At a point in the overlap
of the charts of $A$ and $C$, the representative $A+\mathbf b$ belongs
to $\mathcal U_C$, since that neighborhood is saturated. The same
holds for $\mathbf b$ in an open subset of the domain of $\chi_A$.
On that open set the transition map is
\[
 (\chi_C\circ\chi_A^{-1})(\mathbf b)
   =\operatorname{pr}_{\mathcal S_C}\bigl(\Theta_C^{-1}(A+\mathbf b)\bigr)-C.
\]
Every map on the right is smooth. Interchanging $A$ and $C$ proves
smoothness of the inverse transition. These charts define the asserted
Hilbert manifold structure. In the coordinates $\Theta_A$, the map
$\pi$ is projection onto the first coordinate followed by
$\pi|_{\mathcal S_A}$, hence is a submersion.
At the center of the chart its differential is the projection $\Pi_A$
in \eqref{eq:cotas-proyecciones-coulomb-brecha}, with kernel
$d_AH^{r+1}(E)$ and complement $Y_A$.
\end{proof}

\pagebreak
\begin{proposition}[Full stabilizers and fixed sets]
\label{prop:estabilizadores-completos-fijos-gauge}
Under the hypotheses of
Proposition~\ref{prop:coulomb-sobolev-brecha-uniforme}, assume in addition
that $M$ is connected and noncompact. For
$A\in\mathcal A^r_{A_0}(\mathbf P)$, let
\[
 S(A):=\{\mathbf u\in\mathcal G^{r+1}(A_0)\mid
                                      A\triangleleft\mathbf u=A\}.
\]
Fix $p\in\mathbf P$ and set
$H_A:=Z_G(\operatorname{Hol}_p(A))$. Evaluation induces a group isomorphism
\begin{equation}
\label{eq:estabilizador-completo-holonomia-sobolev}
 \operatorname{ev}_p:S(A)\longrightarrow H_A.
\end{equation}
The group $H_A$ is a compact Lie subgroup of $G$, whereas $S(A)$ is
discrete in the topology induced by $\mathcal G^{r+1}(A_0)$.
The compact topology of $H_A$, transported through
\eqref{eq:estabilizador-completo-holonomia-sobolev}, gives a continuous
action by linear isometries on $H^r(M,F)$ with the norm $\|\cdot\|_{r,A}$.

For $S=S(A)$, the set of connections fixed by all its elements is
the closed affine subspace
\begin{equation}
\label{eq:conjunto-fijo-estabilizador-afine}
 \begin{aligned}
 \mathcal F_S
 &:=\{B\in\mathcal A^r_{A_0}(\mathbf P)\mid S\subseteq S(B)\}
   =A+V_S,\\
 V_S&:=\{\mathbf b\in H^r(M,F)\mid
       \operatorname{Ad}_{\mathbf u^{-1}}\mathbf b=\mathbf b
                         \text{ for every }\mathbf u\in S\}.
 \end{aligned}
\end{equation}
In particular, $\mathcal F_S$ is a split Hilbert submanifold.
\end{proposition}

\begin{proof}
The uniform embedding, with $r>m/2+2$, gives $C^2$ coefficients for $A$.
Parallel transport along piecewise $C^1$ paths and the evaluation
construction in
Theorem~\ref{teo:secciones-paralelas-centralizador-holonomia} apply with
this regularity: the transport equation has coefficients continuous
in the parameter and $C^2$ in the local variables. The construction
and path independence use only uniqueness for that equation and
commutation with holonomy. Thus each $h\in H_A$ determines a unique
$A$-parallel section $\mathbf u_h$ of $\operatorname{Ad}(\mathbf P)$
whose value at $p$ is $h$.

We check its regularity and the global condition for group membership.
On a relatively compact coordinate domain $U$, write the fixed unitary representation as $G\subseteq U(d)$. Then a parallel section satisfies
locally
\[
 \partial_j\mathbf u_h
    =\mathbf u_h A_j-A_j\mathbf u_h,
 \qquad j\in\{1,\ldots,m\}.
\]
The transport construction initially gives a $C^1$ section, hence
an $H^1(U,\mathbb C^{d\times d})$ section on every relatively compact domain.
For $j\in\{1,\ldots,r\}$, the local product
$H^r(U,\mathbb C^{d\times d})\cdot H^j(U,\mathbb C^{d\times d})\subseteq H^j(U,\mathbb C^{d\times d})$ shows that, if
$\mathbf u_h\in H^j_{\mathrm{loc}}(U,\mathbb C^{d\times d})$, all its first derivatives belong
to $H^j_{\mathrm{loc}}(U,\mathbb C^{d\times d})$. It gains one derivative at each step and
therefore belongs to $H^{r+1}_{\mathrm{loc}}(M,\operatorname{Ad}\mathbf P)$.
Cutoffs on smaller domains allow us to apply
Lemma~\ref{lem:modulo-sobolev-uniforme-entero} at each step.

Write $R_h\mathbf b=\operatorname{Ad}_{\mathbf u_h^{-1}}\mathbf b$.
The endomorphism section $R_h$ is $A$-parallel and fiberwise orthogonal.
For $j\in\{0,\ldots,r\}$, in the weak sense,
\[
 (\nabla^A)^j(R_h\mathbf b)=R_h(\nabla^A)^j\mathbf b,
 \qquad \mathbf b\in H^r(M,F),
\]
where $R_h$ acts only on the adjoint factor of each tensor bundle.
The identity is first checked for compactly supported sections using
the local regularity above, then extended by density. It follows that
\begin{equation}
\label{eq:isometria-estabilizador-norma-conexion}
 \|R_h\mathbf b\|_{r,A}=\|\mathbf b\|_{r,A}.
\end{equation}
Theorem~\ref{teo:cambio-conexion-diferencia-sobolev} identifies this
space with $H^r(M,F)$, with equivalent norms. Since
$A=A_0+\mathbf a$ and $d_A\mathbf u_h=0$, the transformation formula gives
\[
 \mathbf u_h^{-1}d_{A_0}\mathbf u_h
       =\mathbf a-R_h\mathbf a\in H^r(M,F).
\]
Thus $\mathbf u_h\in\mathcal G^{r+1}(A_0)$ and stabilizes $A$.
Conversely, the same formula shows that every element of $S(A)$ is
$A$-parallel. Evaluation and extension by transport are inverse maps
and preserve multiplication, proving
\eqref{eq:estabilizador-completo-holonomia-sobolev}.
Compactness of $H_A$ follows from
Proposition~\ref{prop:centralizador-subgrupo-cerrado-lie}.

To check that $S(A)$ is discrete, fix $\iota>0$ such that
$\operatorname{Exp}:B_{\mathfrak g}(0,\iota)\to G$ is a diffeomorphism
onto its image, and an embedding constant $C_{\mathrm{emb}}>0$ for
$H^{r+1}(M,E)\hookrightarrow C_b^0(M,E)$. Choose $\delta>0$ such that
$B_{H^{r+1}(M,E)}(0,\delta)$ lies in the domain of the exponential chart
at the identity and $C_{\mathrm{emb}}\delta<\iota/2$. If
$\mathbf u=\operatorname{Exp}\boldsymbol\xi$ also belongs to $S(A)$,
with $\boldsymbol\xi\in B_{H^{r+1}(M,E)}(0,\delta)$, equivariance of the logarithm
shows that $\boldsymbol\xi$ is $A$-parallel. The infinite volume
argument in
Proposition~\ref{prop:persistencia-espectral-componente-gauge} gives
$\boldsymbol\xi=0$. The identity is therefore isolated in $S(A)$,
and group translations isolate every other element.
In particular, if $H_A$ is infinite, the group isomorphism
\eqref{eq:estabilizador-completo-holonomia-sobolev} is not a
homeomorphism for the Sobolev topology on $S(A)$ and the compact
topology on $H_A$.

Let $H^r(M,F;\nabla^A)$ denote $H^r(M,F)$ equipped with $\|\cdot\|_{r,A}$.
Parallel transport, which is a fiberwise isometry, identifies
$R_h-R_\ell$ at any point with
$\operatorname{Ad}_{h^{-1}}-\operatorname{Ad}_{\ell^{-1}}$ in the
fiber based at $p$. Since this endomorphism is parallel,
\[
 \|R_h-R_\ell\|_{\mathcal L(H^r(M,F;\nabla^A))}
 \leq\|\operatorname{Ad}_{h^{-1}}-\operatorname{Ad}_{\ell^{-1}}\|_{
                                                   \mathcal L(\mathfrak g)},
 \qquad h,\ell\in H_A.
\]
Continuity of the right-hand side proves operator norm continuity
of the action parametrized by the compact group $H_A$.

For $B=A+\mathbf b$ and $\mathbf u\in S$,
$B\triangleleft\mathbf u=A+\operatorname{Ad}_{\mathbf u^{-1}}\mathbf b$.
This proves \eqref{eq:conjunto-fijo-estabilizador-afine}.
The space $V_S$ is closed because it is the intersection, over
$\mathbf u\in S$, of the kernels of the continuous operators
$\operatorname{Ad}_{\mathbf u^{-1}}-I$.
Its complement has an explicit description. If $\mu_A$ is normalized
Haar measure on $H_A$, define
\begin{equation}
\label{eq:proyeccion-promedio-estabilizador-completo}
 P_S\mathbf b:=\int_{H_A}R_h\mathbf b\,d\mu_A(h),
 \qquad\mathbf b\in H^r(M,F).
\end{equation}
The Bochner integral exists by the continuity just proved and
compactness of the domain. Haar invariance gives
$R_\ell P_S\mathbf b=P_S\mathbf b$ for each $\ell\in H_A$; when
$\mathbf b\in V_S$, the integrand is constantly $\mathbf b$.
Hence $P_S^2=P_S$ and $\operatorname{Ran}P_S=V_S$.
Moreover, $R_h^*=R_{h^{-1}}$ for the inner product defining
$\|\cdot\|_{r,A}$. Invariance of Haar measure under inversion gives
$P_S^*=P_S$, and \eqref{eq:isometria-estabilizador-norma-conexion} yields
$\|P_S\|_{\mathcal L(H^r(M,F;\nabla^A))}\leq1$.
If $\Lambda_A\geq1$ satisfies
\[
 \Lambda_A^{-1}\|\mathbf b\|_{H^r(M,F)}
 \leq\|\mathbf b\|_{r,A}\leq\Lambda_A\|\mathbf b\|_{H^r(M,F)},
 \qquad\mathbf b\in H^r(M,F),
\]
then $\|P_S\|_{\mathcal L(H^r(M,F))}\leq\Lambda_A^2$.
The decomposition $H^r(M,F)=V_S\oplus\ker P_S$ proves that the affine
subspace $\mathcal F_S=A+V_S$ is a split submanifold.
\end{proof}

\begin{lemma}[Nearby compact subgroups]
\label{lem:separacion-supergrupos-compactos}
Let $G\subset U(N)$ be a compact Lie group and let $H\subseteq G$ be a
closed subgroup. For the distance
$d_G(g,h)=\|g-h\|_{\mathrm{HS}}$, there is $\eta_H>0$ such that every
closed subgroup $K$ satisfying
\[
 H\subseteq K\subseteq\{g\in G\mid d_G(g,H)<\eta_H\}
\]
equals $H$. Here $\|\cdot\|_{\mathrm{HS}}$ is the Hilbert--Schmidt norm
of matrices and
$\displaystyle d_G(g,H):=\inf_{h\in H}d_G(g,h)$.
\end{lemma}

\begin{proof}
Let $\mathfrak h\subseteq\mathfrak g$ be the Lie algebras. By
Theorem~\ref{teo:subgrupo-cerrado-lie}, we may choose $a>0$ and $L\geq1$
such that the logarithm is defined on
$U:=\{g\in G\mid\|g-I\|_{\mathrm{HS}}<3a\}$, is $L$--Lipschitz there,
$\log(H\cap U)\subseteq\mathfrak h$, and
$\log(\exp X)=X$ for $\|X\|_{\mathrm{HS}}\leq a$.
The Lipschitz bound follows by shrinking a logarithm chart and bounding
its differential on the closure of a slightly larger chart.
Choose $0<\eta_0<a/(2L)$. If
$H\subseteq K\subseteq\{d_G(\cdot,H)<\eta_0\}$ and
$\operatorname{Lie}K$ contains a unit vector
$X\in\mathfrak h^\perp$, then $\exp(aX)\in K$ and
$\|\exp(aX)-I\|_{\mathrm{HS}}\leq a$.
There is $h\in H$ with $\|\exp(aX)-h\|_{\mathrm{HS}}<\eta_0$; hence
$h\in U$ and
\[
 a\leq\|aX-\log h\|_{\mathrm{HS}}
   \leq L\|\exp(aX)-h\|_{\mathrm{HS}}<a/2,
\]
a contradiction. Since $\mathfrak h\subseteq\operatorname{Lie}K$,
the two Lie algebras agree. Their identity components agree as well,
since they are generated by the same exponential neighborhood:
$K_0=H_0$.

Every element of $K$ normalizes $K_0=H_0$. Thus
$K\subseteq N_G(H_0)$, where
$N_G(H_0):=\{g\in G\mid gH_0g^{-1}=H_0\}$ is closed.
The quotient $Q:=N_G(H_0)/H_0$ is a compact Lie group: exponential
charts on a complement of $\mathfrak h$ in the Lie algebra of the
normalizer give quotient charts, and normality of $H_0$ allows the
group operations to descend.
Write $q:N_G(H_0)\to Q$ for the projection. The formula
\[
 d_Q(q(g_1),q(g_2)):=\min_{h\in H_0}\|g_1-g_2h\|_{\mathrm{HS}},
 \qquad g_1,g_2\in N_G(H_0),
\]
defines a bi-invariant distance compatible with the quotient topology.
Compactness of $H_0$ ensures that the minimum is attained and that zero
distance means equality of classes; bi-invariance of $d_G$ and normality
of $H_0$ give independence of representatives and the triangle inequality.
The projection is continuous for this distance; compactness of the
quotient then shows that the two topologies agree.
Set $F=q(H)=H/H_0$, which is finite.

There is $b>0$ such that the ball $B_Q(e,b)$ contains no subgroup other
than $\{e\}$. To see this, choose $\rho>0$ such that the exponential
map of $Q$ is injective on the ball of radius $2\rho$ in its Lie algebra.
No nontrivial subgroup can be contained in $\exp(B(0,\rho))$: if it
contains $\exp X\ne e$ with $0<\|X\|<\rho$, the first integer $n\geq1$
such that $n\|X\|\geq\rho$ satisfies
$\rho\leq n\|X\|<2\rho$, and $\exp(nX)$ does not belong to
$\exp(B(0,\rho))$. It suffices to choose $b>0$ with
$B_Q(e,b)\subseteq\exp(B(0,\rho))$.

If $F$ has more than one element, set
\[
 \mu:=\min_{\substack{f_1,f_2\in F\\f_1\ne f_2}}d_Q(f_1,f_2)>0;
\]
if $F=\{e\}$, set $\mu=1$. Fix
$0<\varepsilon<\min\{b,\mu/4\}$.
For any subgroup $J$ such that
$F\subseteq J\subseteq\{z\in Q\mid d_Q(z,F)<\varepsilon\}$, every
$z\in J$ has distance less than $\varepsilon$ from a unique $f\in F$.
Denote that element by $\tau(z)$.
Bi-invariance gives
\[
 d_Q\bigl(\tau(z_1z_2),\tau(z_1)\tau(z_2)\bigr)<3\varepsilon<\mu,
 \qquad z_1,z_2\in J.
\]
Thus $\tau$ is a homomorphism. Its kernel lies in $B_Q(e,b)$ and is
trivial. Since $\tau|_F=I_F$, for each $z\in J$ the element
$z\tau(z)^{-1}$ belongs to the kernel; hence $J=F$.

Compactness and continuity of $q$ provide $\eta_1>0$ such that
\[
 g\in N_G(H_0),\quad d_G(g,H)<\eta_1
 \quad\Longrightarrow\quad d_Q(q(g),F)<\varepsilon.
\]
Indeed, the closed set where the last distance is at least
$\varepsilon$ is compact and disjoint from $H$, so its distance from
$H$ is positive; if it is empty, any $\eta_1>0$ works.
Taking $\eta_H=\min\{\eta_0,\eta_1\}$, the preceding steps give
$q(K)=F$ and $K_0=H_0$, and therefore $K=H$.
\end{proof}

\begin{proposition}[Openness and density of the exact stabilizer set]
\label{prop:estabilizador-exacto-abierto-denso}
Under the hypotheses and notation of
Proposition~\ref{prop:estabilizadores-completos-fijos-gauge}, fix
$A\in\mathcal A^r_{A_0}(\mathbf P)$ and $S=S(A)$. Then
\[
 \mathcal C_S:=\{B\in\mathcal A^r_{A_0}(\mathbf P)\mid S(B)=S\}
\]
is open and dense in $\mathcal F_S=A+V_S$. In particular,
$\mathcal C_S$ is a split Hilbert submanifold and
\begin{equation}
\label{eq:clausura-estabilizador-exacto-fijo}
 \overline{\mathcal C_S}^{\,\mathcal A^r_{A_0}(\mathbf P)}=\mathcal F_S.
\end{equation}
\end{proposition}

\begin{proof}
Fix $x_*\in M$ and $p\in\mathbf P_{x_*}$, and set
$H=\operatorname{ev}_p(S)=Z_G(\operatorname{Hol}_p(A))$.
For $B\in\mathcal F_S$, evaluation identifies $S(B)$ with the compact
subgroup $H_B:=Z_G(\operatorname{Hol}_p(B))$, which contains $H$.
Moreover, $S(B)=S$ if and only if $H_B=H$: sections in both groups are
$B$-parallel, and a parallel section is determined by its value at $p$.
If $H=G$, this equality holds for every $B\in\mathcal F_S$; henceforth
we may assume $H\ne G$.

We write $\operatorname{Hol}_p(B;\gamma)$ for the element of $G$
determined by the transport of $B$ along the loop $\gamma$.
First we prove openness at a point $B\in\mathcal C_S$.
There are piecewise smooth loops $\gamma_1,\ldots,\gamma_q$ based at
$x_*$ such that, writing $g_i=\operatorname{Hol}_p(B;\gamma_i)$,
\[
 H=\bigcap_{i=1}^q Z_G(g_i).
\]
Indeed, the commutant of the holonomy in
$\operatorname{End}(\mathbb C^N)$ is an intersection of kernels of
linear operators on a finite-dimensional space. Successively choose
a condition that reduces the dimension of the partial intersection;
after at most $N^2$ reductions, we obtain the full commutant.
Intersecting with $G$ gives the displayed equality.

Let $\eta_H$ be the constant in
Lemma~\ref{lem:separacion-supergrupos-compactos}. If the compact set
$D:=\{g\in G\mid d_G(g,H)\geq\eta_H\}$ is nonempty, define
\[
 \alpha:=\min_{g\in D}\max_{1\leq i\leq q}
                         \|gg_i-g_i g\|_{\mathrm{HS}}>0.
\]
For each loop, the transport equation and unitarity give
\begin{equation}
\label{eq:control-holonomia-perturbacion-sobolev}
 \|\operatorname{Hol}_p(C;\gamma_i)-g_i\|_{\mathrm{HS}}
 \leq c_{\mathrm{rep}}\ell_{\mathbf g}(\gamma_i)
                         \|C-B\|_{C_b^0(M,F)}
 \leq L_i\|C-B\|_{H^r(M,F)}.
\end{equation}
Here $c_{\mathrm{rep}}>0$ bounds the differential of the unitary
representation, $C_{\mathrm{emb}}>0$ is the constant of
$H^r(M,F)\hookrightarrow C_b^0(M,F)$, and
$L_i=c_{\mathrm{rep}}\ell_{\mathbf g}(\gamma_i)C_{\mathrm{emb}}$.
In a unitary trivialization along the loop, write $U_C(t,s)$ for
transport from $s$ to $t$, with $0\leq s\leq t\leq1$.
The first inequality follows by expressing the difference of the
transports as
\[
 -\int_0^1 U_C(1,t)(C-B)_{\gamma_i(t)}(\dot\gamma_i(t))
                                  U_B(t,0)\,dt.
\]
The factors $U_C$ and $U_B$ have operator norm one. The identity applies
on each piece, and the contributions sum to the total length.

Set $L=\max\{1,L_1,\ldots,L_q\}$. If $C\in\mathcal F_S$ and
$\|C-B\|_{H^r(M,F)}<\alpha/(4L)$, every $g\in H_C$ satisfies
\[
 \max_{1\leq i\leq q}\|gg_i-g_i g\|_{\mathrm{HS}}<\alpha/2,
\]
since it commutes with $\operatorname{Hol}_p(C;\gamma_i)$ and is
unitary. Thus
$H\subseteq H_C\subseteq\{d_G(\cdot,H)<\eta_H\}$.
The lemma gives $H_C=H$, and hence $C\in\mathcal C_S$.
If $D$ is empty, the lemma gives the same conclusion for any
$C\in\mathcal F_S$, without restricting its distance from $B$.
This proves relative openness.

To prove density, now fix any $B\in\mathcal F_S$ and $\varepsilon>0$.
Write $C_H:=Z_G(H)$ and $\mathfrak c:=\operatorname{Lie}C_H$.
Every holonomy transport of a connection in $\mathcal F_S$ belongs
to $C_H$, since its parallel sections include $S$.
For each loop $\gamma$, the affine segment between two connections
in $\mathcal F_S$ and
\eqref{eq:control-holonomia-perturbacion-sobolev} show that their
holonomies along $\gamma$ lie in the same connected component of $C_H$.
If $\mathfrak c=\{0\}$, the group $C_H$ is finite; these holonomies
equal those of $A$ for every loop, and $H_B=H$.
If $m=1$, the connected noncompact base without boundary is
diffeomorphic to $\mathbb R$; every connection has trivial holonomy,
and $H=G$, a case already considered. It remains to treat
$d:=\dim\mathfrak c>0$ and $m\geq2$.

We construct a local trivialization in which every element of $S$ is
constant. Choose $s\geq1$ and $h_1,\ldots,h_s\in H$ such that
$\displaystyle C_H=\bigcap_{j=1}^s Z_G(h_j)$; the existence of this
finite family follows from the same commutant argument used above,
adding the identity if the initial family is empty.
In a smooth trivialization around $x_*$ whose value at $x_*$ is $p$,
the sections $\mathbf u_{h_j}\in S$ have representatives $U_j(x)$ of
class $H^{r+1}_{\mathrm{loc}}(U,G)$.
The tuple $(U_1(x),\ldots,U_s(x))$ lies in the conjugation orbit of
$(h_1,\ldots,h_s)$, since these sections are $B$-parallel.
The map
\[
 \Phi:G\longrightarrow G^s,\qquad
 \Phi(g)=(g h_1g^{-1},\ldots,g h_sg^{-1})
\]
has differential kernel $\mathfrak c$ at the identity and constant
rank by equivariance. Theorem~\ref{teo: teo del rango}, applied on a complement of
$\mathfrak c$, gives a local inverse onto its image.
Compactness of $G$ ensures that this image contains a neighborhood of
$\Phi(I)$ in the orbit: if $\Phi(g_n)\to\Phi(I)$, a subsequence
satisfies $g_n\to c\in C_H$, and $\Phi(g_n c^{-1})=\Phi(g_n)$ with
$g_n c^{-1}\to I$.
Consequently, after shrinking the neighborhood of $x_*$, there is
$g(x)\in H^{r+1}_{\mathrm{loc}}(U,G)$ with $g(x_*)=I$ and
$U_j(x)=g(x)h_jg(x)^{-1}$ for $j=1,\ldots,s$.
Regularity follows from
Lemma~\ref{lem:superposicion-sobolev-geometrica}, applied locally to
the smooth inverse of $\Phi$ on that image.

Changing the trivialization by $g(x)$ makes every element of $S$
constant: two elements of $G$ that conjugate the chosen tuple in the
same way differ by an element of $C_H$, which commutes with all of $H$.
Since $S$ is $B$-parallel, the local coefficients of $B$ in this
trivialization take values in $\mathfrak c$.
These coefficients belong to $H^r_{\mathrm{loc}}(U,T^*U\otimes\mathfrak c)$ by the connection
transformation formula and the Sobolev product estimates.

Choose coordinates $x=(x^1,\ldots,x^m)$ on this neighborhood with
$x_*=0$, and a basis $X_1,\ldots,X_d$ of $\mathfrak c$.
Choose $R>0$ such that the closed Euclidean ball of radius $R$ lies
in the coordinate domain, and set
\[
 t_i:=\frac{iR}{2(d+1)},\qquad q_i:=(t_i,0,\ldots,0),\qquad
 \delta:=\frac{R}{8(d+1)},\qquad D_i:=B_{\mathbb R^m}(q_i,\delta),
 \quad i=1,\ldots,d.
\]
The balls $D_i$ have disjoint closures, lie inside the chart, and do
not contain the origin. Let $\chi_i\in C_c^\infty(D_i)$ equal one on
$B_{\mathbb R^m}(q_i,\delta/2)$, for $i=1,\ldots,d$.
Let $T_i\in C_H$ be the transport of $B$ along the segment
$\lambda_i(s)=(s t_i,0,\ldots,0)$, $0\leq s\leq1$, from the origin to
$q_i$, written in the chosen trivialization, and
set $Y_i=\operatorname{Ad}_{T_i}X_i$.
Within the chart define
\begin{equation}
\label{eq:perturbacion-curvatura-estabilizador-exacto}
 \mathbf b:=\sum_{i=1}^d
       \chi_i(x)(x^1-t_i)Y_i\,dx^2,
\end{equation}
and extend by zero outside it. This section belongs to $H^r(M,F)$:
its expression in the preceding trivialization is smooth and compactly
supported, and the change to a smooth trivialization is a multiplier
of class $H^{r+1}_{\mathrm{loc}}(U,G)$.
Moreover, $\mathbf b\in V_S$, since $Y_i\in\mathfrak c$.

For $B_t:=B+t\mathbf b$, with $t\in\mathbb R$, all transports along
the segments $\lambda_i$ remain equal to $T_i$: we have
$\mathbf b(\dot\lambda_i)=0$ because $dx^2(\dot\lambda_i)=0$.
At $q_i$, the section $\mathbf b$ vanishes and
$(d_B\mathbf b)_{q_i}(\partial_1,\partial_2)=Y_i$.
The curvature formula therefore gives the exact identity
\[
 \operatorname{Ad}_{T_i^{-1}}
       \mathbf F_{B_t,q_i}(\partial_1,\partial_2)
 =Z_i+tX_i,\qquad
 Z_i:=\operatorname{Ad}_{T_i^{-1}}
       \mathbf F_{B,q_i}(\partial_1,\partial_2)\in\mathfrak c,
 \quad i=1,\ldots,d.
\]
In the basis $(X_1,\ldots,X_d)$, the determinant of the vectors
$(Z_1+tX_1,\ldots,Z_d+tX_d)$ is a polynomial of degree $d$ with leading
coefficient one. It therefore has finitely many zeros. We may choose
\[
 0<t<\frac{\varepsilon}{1+\|\mathbf b\|_{H^r(M,F)}}
\]
outside those zeros. For this choice the preceding vectors span
$\mathfrak c$, and $\|B_t-B\|_{H^r(M,F)}<\varepsilon$.

If $\mathbf u\in S(B_t)$ and $z=\operatorname{ev}_p(\mathbf u)$, the
identity $d_{B_t}\mathbf u=0$ implies
$[\mathbf F_{B_t},\mathbf u]=0$. Transporting from $q_i$ to $x_*$ shows
that $z$ commutes with each $Z_i+tX_i$.
It therefore commutes with $\mathfrak c$ and with the connected group
$(C_H)_0$, generated by its exponentials.
For any loop $\gamma$ based at $x_*$, the elements
$\operatorname{Hol}_p(A;\gamma)$ and
$\operatorname{Hol}_p(B_t;\gamma)$ lie in the same component of $C_H$,
as proved above. Their quotient belongs to $(C_H)_0$.
The element $z$ commutes with this quotient and with the holonomy of
$B_t$, so it commutes with the holonomy of $A$ along every loop.
Thus $z\in H$ and $H_{B_t}=H$, since the reverse inclusion follows
from $B_t\in\mathcal F_S$. We have constructed $B_t\in\mathcal C_S$
at distance less than $\varepsilon$ from $B$.

This proves density. Since $\mathcal F_S$ is closed in the affine
component, we obtain \eqref{eq:clausura-estabilizador-exacto-fijo}.
Relative openness and the split structure of $\mathcal F_S$ established
in Proposition~\ref{prop:estabilizadores-completos-fijos-gauge} give
the submanifold assertion.
\end{proof}

\begin{corollary}[Density after saturation]
\label{cor:densidad-saturacion-estabilizador-exacto}
Under the hypotheses of the preceding proposition, define
\[
 \begin{aligned}
 \mathcal C_{(S)}&:=\{B\in\mathcal A^r_{A_0}(\mathbf P)\mid S(B)\in(S)\},\\
 \mathcal D_S&:=\bigcup_{\mathbf u\in\mathcal G_\bullet^{r+1}(A_0)}
                                      (\mathcal F_S\triangleleft\mathbf u),
 \end{aligned}
\]
where $(S)$ is the
conjugacy class of $S$ under $\mathcal G_\bullet^{r+1}(A_0)$.
Then $\mathcal C_{(S)}$ is dense in $\mathcal D_S$ for the relative
topology, and
\begin{equation}
\label{eq:clausuras-saturaciones-estabilizadores}
 \overline{\mathcal C_{(S)}}=\overline{\mathcal D_S}.
\end{equation}
Both closures are taken in the affine component.
\end{corollary}

\begin{proof}
The identity $S(B\triangleleft\mathbf u)=\mathbf u^{-1}S(B)\mathbf u$
gives
\[
 \mathcal C_{(S)}
 =\bigcup_{\mathbf u\in\mathcal G_\bullet^{r+1}(A_0)}
                              (\mathcal C_S\triangleleft\mathbf u)
 \subseteq\mathcal D_S.
\]
Every translation is a homeomorphism. By the preceding proposition,
$\mathcal F_S\triangleleft\mathbf u$ is contained in the closure of
$\mathcal C_S\triangleleft\mathbf u$ for every $\mathbf u$.
Taking the union gives
$\mathcal D_S\subseteq\overline{\mathcal C_{(S)}}$, proving both claims.
\end{proof}

\pagebreak
\begin{lemma}[Comparison of parallel families]
\label{lem:comparacion-familias-paralelas-gauge}
Assume the hypotheses of
Corollary~\ref{cor:cociente-hilbert-brecha-gauge}, and fix a connection
$A_*\in\mathcal A^r_{A_0}(\mathbf P)$ and $p\in\mathbf P_{x_*}$.
Let $s\geq1$, $h_1,\ldots,h_s\in G$, and
\[
 \mathcal Q:=\{(g h_1g^{-1},\ldots,g h_sg^{-1})\mid g\in G\}.
\]
Consider the associated bundle $\mathbf P\times_G\mathcal Q$ for
the diagonal conjugation action. There are a neighborhood
$\mathcal V_*$ of $A_*$ and constants $\delta>0$, $a>0$, and $J\geq1$
with the following property. Let $A,B\in\mathcal V_*$, and let
$\boldsymbol\sigma$, $\boldsymbol\tau$ be sections of this bundle
of class $H^{r+1}_{\mathrm{loc}}(M,\mathbf P\times_G\mathcal Q)$, parallel with respect to $A$ and
$B$, respectively. If
\[
 \boldsymbol\tau-\boldsymbol\sigma\in
       L^2\bigl(M,(\operatorname{End}\mathbf E)^s\bigr),\qquad
 \|B-A\|_{H^r(M,F)}<\delta,
 \qquad d_{\mathcal Q}(\boldsymbol\sigma(p),\boldsymbol\tau(p))<a,
\]
then there is $\boldsymbol\xi\in H^{r+1}(M,E)$ such that
\begin{equation}
\label{eq:comparacion-familias-paralelas-gauge}
 \boldsymbol\tau
   =\operatorname{Exp}(\boldsymbol\xi)\boldsymbol\sigma
                         \operatorname{Exp}(-\boldsymbol\xi),
 \qquad
 \|\boldsymbol\xi\|_{H^{r+1}(M,E)}
       \leq J\|B-A\|_{H^r(M,F)}.
\end{equation}
Conjugation applies to each member of the family.
Here $\mathbf E$ is the bundle of the fixed unitary representation,
and $d_{\mathcal Q}$ is an invariant distance which, on each connected
component of $\mathcal Q$, comes from the homogeneous metric described
in the proof. Distinct components are assigned a constant distance
greater than all their diameters.
\end{lemma}

\begin{proof}
The stabilizer of the tuple $(h_1,\ldots,h_s)$ is
$\displaystyle C=\bigcap_{j=1}^s Z_G(h_j)$.
The rank theorem and compactness of $G$ smoothly identify
$\mathcal Q$ with $G/C$, as in the trivialization construction of
Proposition~\ref{prop:estabilizador-exacto-abierto-denso}.
Fix an $\operatorname{Ad}(G)$--invariant inner product on $\mathfrak g$, and give the
quotient its homogeneous metric. For $q\in\mathcal Q$, write
$\mathfrak c_q$ for the Lie algebra of its stabilizer. On
$\mathfrak c_q^\perp$, the map
\[
 X\longmapsto\exp(X)q\exp(-X)
\]
is a normal chart at $q$: the horizontal group geodesic $\exp(tX)$
projects to the quotient geodesic. Homogeneity and compactness provide
$\iota>0$ that works for every point of $\mathcal Q$.
For $d_{\mathcal Q}(q,q')<\iota$, write
$L(q,q')\in\mathfrak c_q^\perp$ for the unique vector in this chart.
Then
\[
 |L(q,q')|=d_{\mathcal Q}(q,q'),\qquad
 q'=\exp(L(q,q'))q\exp(-L(q,q')),
\]
and $L$ is smooth and equivariant under conjugation. If
$\mathcal Q$ has dimension zero, take $\iota$ smaller than the
separation between distinct points; the same construction applies
with $L(q,q)=0$.

Choose $0<\eta<\iota/4$ and a smooth function
$\chi:[0,\infty)\to[0,1]$ equal to one on $[0,4\eta^2]$ and zero
on $[9\eta^2,\infty)$. The map
\[
 f(q,q'):=\chi(|L(q,q')|^2)L(q,q')
\]
extends by zero to $\mathcal Q\times\mathcal Q$ and is smooth and
equivariant. Its derivatives are bounded and $f(q,q)=0$.
In the matrix representation, there is $C_f>0$ such that
\[
 |f(q,q')|\leq C_f|q'-q|,
 \qquad q,q'\in\mathcal Q,
\]
using the Euclidean norm of the Hilbert--Schmidt norms for tuples.
The bound follows near the diagonal from the differential of $f$,
and away from that neighborhood from compactness.

We may choose $\mathcal V_*$ to be a ball of radius smaller than
$\min\{1,\sqrt{\gamma_{A_*}}/(2c)\}$, with $c$ from
Proposition~\ref{prop:persistencia-espectral-componente-gauge}.
The gap is at least $\gamma_{A_*}/4$ throughout this ball.
The bounds in \eqref{eq:cotas-proyecciones-coulomb-brecha} and uniform
equivalence of norms for connections in a bounded ball provide
$C_*\geq1$ such that
\begin{equation}
\label{eq:coercividad-uniforme-familias-paralelas}
 \|\boldsymbol\zeta\|_{H^{r+1}(M,E)}
 \leq C_*\|d_A\boldsymbol\zeta\|_{H^r(M,F)},
 \qquad A\in\mathcal V_*,\quad\boldsymbol\zeta\in H^{r+1}(M,E).
\end{equation}

Set $\boldsymbol\beta=B-A$ and
$\boldsymbol\xi=f(\boldsymbol\sigma,\boldsymbol\tau)$.
Equivariance makes $\boldsymbol\xi$ a section of $E$.
The integrability assumption and the bound for $f$ give
$\boldsymbol\xi\in L^2(M,E)$. Moreover,
\[
 d_A\boldsymbol\sigma=0,\qquad
 d_A\boldsymbol\tau=-[\boldsymbol\beta,\boldsymbol\tau],\qquad
 d_A\boldsymbol\xi
   =D_2f(\boldsymbol\sigma,\boldsymbol\tau)
                         \bigl(-[\boldsymbol\beta,\boldsymbol\tau]\bigr).
\]
Each additional derivative either differentiates a factor of
$\boldsymbol\beta$ or adds another when it differentiates a coefficient
depending on $\boldsymbol\tau$; coefficients depending on
$\boldsymbol\sigma$ are parallel. The induction in
\eqref{eq:logaritmo-truncado-calibre}, applied to this two-variable
map, gives, for $j=1,\ldots,r+1$, constants $D_j>0$ such that
\[
 |(\nabla^A)^j\boldsymbol\xi|
 \leq D_j\sum_{q=1}^{j}
 \sum_{\substack{\alpha_1,\ldots,\alpha_q\geq0\\
                  \alpha_1+\cdots+\alpha_q=j-q}}
 \prod_{\nu=1}^{q}|(\nabla^A)^{\alpha_\nu}\boldsymbol\beta|.
\]
The indices in the inner sum are integers. The local chain rule
and Sobolev approximation justify these identities for weak
derivatives. The bounds on derivatives of $f$ are uniform on the
compact set $\mathcal Q\times\mathcal Q$. The Sobolev module estimates
and equivalence of norms, uniform for $A\in\mathcal V_*$, provide
$C\geq1$ such that
\[
 \|d_A\boldsymbol\xi\|_{H^r(M,F)}
       \leq C\sum_{q=1}^{r+1}\|\boldsymbol\beta\|_{H^r(M,F)}^q.
\]
Together with the order-zero term already controlled, these estimates
prove $\boldsymbol\xi\in H^{r+1}(M,E)$ before applying
\eqref{eq:coercividad-uniforme-familias-paralelas}.

Let $C_{\mathrm{emb}}>0$ be the constant of
$H^{r+1}(M,E)\hookrightarrow C_b^0(M,E)$, and set
\[
 J:=C_*C(r+1),\qquad
 \delta:=\min\{1/2,\eta/(2C_{\mathrm{emb}}J)\},\qquad a:=\eta/2.
\]
Under the hypotheses, we obtain
$\|\boldsymbol\xi\|_{H^{r+1}(M,E)}\leq J\|\boldsymbol\beta\|_{H^r(M,F)}$
and $\|\boldsymbol\xi\|_\infty<\eta/2$.
The distance between the two families is independent of the
trivialization, continuous on $M$, and less than $\eta/2$ at $x_*$.
If it reached $\eta$ elsewhere, a path from $x_*$ would have a first
point where the distance equals $\eta$. There $|f|=\eta$, a
contradiction. Thus the distance is less than $\eta$ throughout the
base, and $f=L$ on the two families. This proves
\eqref{eq:comparacion-familias-paralelas-gauge}.
The exponential of $\boldsymbol\xi$ belongs to
$\mathcal G_e^{r+1}(A_0)$, since
$t\mapsto\operatorname{Exp}(t\boldsymbol\xi)$, $0\leq t\leq1$,
is a continuous path of gauge transformations starting at the
identity.
\end{proof}

\begin{proposition}[Closedness of saturated sets]
\label{prop:clausura-saturaciones-estabilizadores-gauge}
Under the hypotheses of
Corollary~\ref{cor:cociente-hilbert-brecha-gauge}, let $S=S(A_1)$
for some $A_1\in\mathcal A^r_{A_0}(\mathbf P)$, and use the notation
of Corollary~\ref{cor:densidad-saturacion-estabilizador-exacto}.
Then $\mathcal D_S$ is closed in the affine component and
\begin{equation}
\label{eq:clausura-tipo-saturacion-gauge}
 \overline{\mathcal C_{(S)}}=\mathcal D_S.
\end{equation}
\end{proposition}

\begin{proof}
Fix $p\in\mathbf P_{x_*}$ and $H=\operatorname{ev}_p(S)$.
Choose $s\geq1$ and $h_1,\ldots,h_s\in H$ such that
$\displaystyle Z_G(H)=\bigcap_{j=1}^s Z_G(h_j)$, as in
Proposition~\ref{prop:estabilizador-exacto-abierto-denso}.
Let $\boldsymbol\sigma=(\mathbf u_{h_1},\ldots,\mathbf u_{h_s})$
be the corresponding $A_1$-parallel sections.
These sections determine the fixed set:
\begin{equation}
\label{eq:familia-finita-determina-fijo-gauge}
 B\in\mathcal F_S
 \quad\Longleftrightarrow\quad
 d_B\mathbf u_{h_j}=0\quad\text{for }j=1,\ldots,s.
\end{equation}
Indeed, writing $B=A_1+\mathbf b$, the last equations say that
$\mathbf b$ commutes with each $\mathbf u_{h_j}$.
In every fiber, the chosen equality of centralizers implies that
it commutes with every element of $S$. Hence
$d_B\mathbf u_h=0$ for all $h\in H$, which is precisely
$S\subseteq S(B)$. The same equivalence holds after any gauge
transformation.

Let $B_n\in\mathcal D_S$, $n\in\mathbb N$, be a sequence converging
to $B$ in $H^r(M,T^*M\otimes\operatorname{ad}\mathbf P)$. The central components act trivially; we may choose
$\mathbf u_n\in\mathcal G_e^{r+1}(A_0)$ such that
$\mathbf u_n^{-1}S\mathbf u_n\subseteq S(B_n)$.
The families
\[
 \boldsymbol\sigma_n:=\mathbf u_n^{-1}\boldsymbol\sigma\mathbf u_n
\]
are $B_n$-parallel and take values in the bundle associated with the
compact orbit $\mathcal Q$ of the tuple $(h_1,\ldots,h_s)$.
Moreover,
\[
 \boldsymbol\sigma_n-\boldsymbol\sigma_m\in
              L^2\bigl(M,(\operatorname{End}\mathbf E)^s\bigr),
 \qquad n,m\in\mathbb N.
\]
To verify this, Lemma~\ref{lem:control-global-calibre-brecha} gives
$\mathbf u_n-I\in H^{r+1}(M,\operatorname{End}\mathbf E)$.
The identity
\[
 \mathbf u_n^{-1}\boldsymbol\sigma\mathbf u_n-\boldsymbol\sigma
 =\mathbf u_n^{-1}\boldsymbol\sigma(\mathbf u_n-I)
       +(\mathbf u_n^{-1}-I)\boldsymbol\sigma
\]
and pointwise boundedness of the matrices prove integrability.

Passing to a subsequence, we may assume that
$\boldsymbol\sigma_n(p)\to q\in\mathcal Q$.
Continuity of transport with respect to the connection, proved in
\eqref{eq:control-holonomia-perturbacion-sobolev}, implies that each
component of $q$ commutes with the entire holonomy of $B$.
By Theorem~\ref{teo:secciones-paralelas-centralizador-holonomia},
with the Sobolev regularity justified in
Proposition~\ref{prop:estabilizadores-completos-fijos-gauge}, $q$
determines a $B$-parallel family $\boldsymbol\tau$.
For every $x\in M$, choosing a fixed path from $x_*$ to $x$ and
applying the same continuity of transport gives
$\boldsymbol\sigma_n(x)\to\boldsymbol\tau(x)$.

Apply Lemma~\ref{lem:comparacion-familias-paralelas-gauge}
with $A_*=B$. For all sufficiently large $n,m$,
$B_n,B_m\in\mathcal V_*$,
$\|B_n-B_m\|_{H^r(M,F)}<\delta$, and
$d_{\mathcal Q}(\boldsymbol\sigma_m(p),\boldsymbol\sigma_n(p))<a$.
Fix one such index $m$. For every sufficiently large $n$, the lemma
provides
\[
 \boldsymbol\xi_n=f(\boldsymbol\sigma_m,\boldsymbol\sigma_n),\qquad
 \boldsymbol\sigma_n
  =\operatorname{Exp}(\boldsymbol\xi_n)\boldsymbol\sigma_m
                          \operatorname{Exp}(-\boldsymbol\xi_n),
 \qquad
 \|\boldsymbol\xi_n\|_{H^{r+1}(M,E)}
       \leq J\|B_n-B_m\|_{H^r(M,F)}.
\]
This sequence is bounded in the Hilbert space $H^{r+1}(M,E)$.
A subsequence converges weakly. Its limit is identified with
$\boldsymbol\xi=f(\boldsymbol\sigma_m,\boldsymbol\tau)$:
the preceding pointwise convergence and the uniform Sobolev embedding
bound give strong $L^2(K,E\restriction_K)$ convergence for every compact set $K\subseteq M$, by dominated
convergence. Therefore
\[
 \boldsymbol\xi\in H^{r+1}(M,E),\qquad
 \|\boldsymbol\xi\|_{H^{r+1}(M,E)}
       \leq J\|B-B_m\|_{H^r(M,F)},\qquad
 \boldsymbol\tau
  =\operatorname{Exp}(\boldsymbol\xi)\boldsymbol\sigma_m
                          \operatorname{Exp}(-\boldsymbol\xi).
\]
Set $\mathbf v=\operatorname{Exp}(\boldsymbol\xi)$ and
$\mathbf w=\mathbf u_m\mathbf v^{-1}\in\mathcal G_e^{r+1}(A_0)$.
Then $\boldsymbol\tau=\mathbf w^{-1}\boldsymbol\sigma\mathbf w$.
Since its components are $B$-parallel, the transformed version of
\eqref{eq:familia-finita-determina-fijo-gauge} gives
$\mathbf w^{-1}S\mathbf w\subseteq S(B)$.
Thus $B\in\mathcal D_S$. The affine component is metrizable, so
this verification with sequences proves that $\mathcal D_S$ is closed.
Equality \eqref{eq:clausura-tipo-saturacion-gauge} follows from
Corollary~\ref{cor:densidad-saturacion-estabilizador-exacto}.
\end{proof}

\begin{proposition}[Structure of the sets of fixed type]
\label{prop:estratos-hilbert-estabilizadores-gauge}
Under the hypotheses of the preceding proposition,
$\mathcal C_{(S)}$ is a split Hilbert submanifold, open and dense in
$\mathcal D_S$, and locally closed in the affine component.
Its image $\mathcal R_{(S)}:=\pi(\mathcal C_{(S)})$ is a split,
locally closed Hilbert submanifold of $\mathcal R$. Moreover,
\begin{equation}
\label{eq:clausura-estrato-cociente-gauge}
 \overline{\mathcal R_{(S)}}=\pi(\mathcal D_S),
\end{equation}
and the projection
$\pi:\mathcal C_{(S)}\to\mathcal R_{(S)}$ is a smooth submersion.
The smooth structure on $\mathcal R_{(S)}$ is the unique one with
this last property and with the relative topology of the quotient.
\end{proposition}

\begin{proof}
Fix $A\in\mathcal C_{(S)}$. We may use $S(A)$ as the representative
of the type, since this replacement changes neither
$\mathcal C_{(S)}$ nor $\mathcal D_S$; write $S=S(A)$ for the local
argument. Take the tube
$\Theta_A:\mathcal S_A\times\mathcal G_e^{r+1}(A_0)\to\mathcal U_A$
from Proposition~\ref{prop:tubo-calibre-brecha-conexa}.
The identity component is normal in the full group, and the slice
is invariant under $S$. Consequently,
\begin{equation}
\label{eq:equivarianza-estabilizador-tubo-gauge}
 \Theta_A(C,\mathbf u)\triangleleft\mathbf h
 =\Theta_A(C\triangleleft\mathbf h,\mathbf h^{-1}\mathbf u\mathbf h),
 \quad C\in\mathcal S_A,\quad
 \mathbf u\in\mathcal G_e^{r+1}(A_0),\quad\mathbf h\in S.
\end{equation}
In particular, $\mathcal U_A$ is invariant under $S$. If
$D\in\mathcal U_A$ is fixed by $S$, uniqueness of the coordinates
$\Theta_A^{-1}(D)=(C,\mathbf u)$ and the preceding identity show
that $C$ is also fixed by $S$.

We obtain the equality
\begin{equation}
\label{eq:saturacion-en-tubo-gauge}
 \mathcal D_S\cap\mathcal U_A
   =\Theta_A\bigl((\mathcal F_S\cap\mathcal S_A)
                             \times\mathcal G_e^{r+1}(A_0)\bigr).
\end{equation}
Indeed, if $D\in\mathcal D_S\cap\mathcal U_A$, some transformation
in the identity component takes $D$ to a point fixed by $S$.
That point remains in $\mathcal U_A$ and has the same representative
$C\in\mathcal S_A$ as $D$. The preceding paragraph gives
$C\in\mathcal F_S$. The reverse inclusion follows from the definition
of $\mathcal D_S$.

Shrink the slice radius $\varepsilon_A$ until
$\mathcal F_S\cap\mathcal S_A\subseteq\mathcal C_S$, using openness
from Proposition~\ref{prop:estabilizador-exacto-abierto-denso}.
Shrinking preserves the tube and its invariance. By
\eqref{eq:saturacion-en-tubo-gauge}, we obtain
\begin{equation}
\label{eq:modelo-local-tipo-estabilizador-gauge}
 \mathcal C_{(S)}\cap\mathcal U_A
 =\mathcal D_S\cap\mathcal U_A
 =\Theta_A\bigl((A+B_{Y_A\cap V_S}(0,\varepsilon_A))
                           \times\mathcal G_e^{r+1}(A_0)\bigr),
\end{equation}
where the ball uses the norm $\|\cdot\|_{r,A}$.

The space $Y_A\cap V_S$ is a closed complemented subspace of $Y_A$.
Indeed, each operator $R_h$ of
Proposition~\ref{prop:estabilizadores-completos-fijos-gauge} preserves
$Y_A=\ker d_A^*$, by the identity $d_A^*R_h=R_h d_A^*$.
The Haar projection $P_S$ restricted to $Y_A$ is a continuous
projection onto $Y_A\cap V_S$.
Formula \eqref{eq:modelo-local-tipo-estabilizador-gauge} proves the
split submanifold assertion and openness in $\mathcal D_S$.
Density has already been proved; closedness of $\mathcal D_S$ shows
that $\mathcal C_{(S)}$ is locally closed.

In the quotient charts constructed in
Corollary~\ref{cor:cociente-hilbert-brecha-gauge}, the image of
\eqref{eq:modelo-local-tipo-estabilizador-gauge} is exactly
$B_{Y_A\cap V_S}(0,\varepsilon_A)$. The chart transitions restrict
the smooth transitions of the ambient quotient, and their inverses
are the corresponding restrictions. They therefore give
$\mathcal R_{(S)}$ its split submanifold structure.
In tube coordinates, the projection onto the stratum is projection
onto this first factor, so it is a smooth submersion.
If another smooth structure with the same topology made $\pi$ a
submersion, its local sections, supplied by the submersion theorem,
would show that the identity between the two structures is smooth
in both directions. This proves uniqueness.

The set $\mathcal D_S$ is closed and saturated; by definition of the
quotient topology, $\pi(\mathcal D_S)$ is closed in $\mathcal R$.
It contains $\mathcal R_{(S)}$ and, by density of
$\mathcal C_{(S)}$ in $\mathcal D_S$ and continuity of $\pi$, it lies
in $\overline{\mathcal R_{(S)}}$. This proves
\eqref{eq:clausura-estrato-cociente-gauge}.
The local equality \eqref{eq:modelo-local-tipo-estabilizador-gauge}
also shows that $\mathcal R_{(S)}$ is open in this closed set,
and is therefore locally closed.
\end{proof}

\pagebreak
\begin{theorem}[Eichhorn--Heber: gauge slice in bounded geometry without boundary]
\label{teo:eichhorn-heber-slice-gauge}
Let $(M,\mathbf{g})$ be a smooth connected noncompact Riemannian manifold of dimension $m$
without boundary and of bounded geometry of order $k$, in the sense of
Definition~\ref{def:geometria-acotada-orden-k}. Let
$\mathbf{P}\longrightarrow M$ be a
principal bundle with compact group $G$, let $A_0$ be a smooth connection satisfying
\eqref{eq:geometria-acotada-conexion-gauge} and
\[
 \operatorname{YM}(A_0)=\frac12\int_M|\mathbf{F}_{A_0}|_{\mathbf{g},h_{\mathrm{ad}}}^2\,d\lambda_{\mathbf{g}}<\infty,
\]
and suppose
\[
k-1\geq r>\frac{m}{2}+2,\qquad k,r\in\mathbb N.
\]
Assume that the reference connection satisfies
\[
 \inf\sigma_{\mathrm e}(\Delta_{A_0}^{(0)})>0.
\]
Then, for every $A\in\mathcal A^r_{A_0}(\mathbf P)$, the operator
$d_A:H^{r+1}(M,\operatorname{ad}\mathbf P)\to
H^r(M,T^*M\otimes\operatorname{ad}\mathbf P)$ has closed range in
$H^r(M,T^*M\otimes\operatorname{ad}\mathbf P)$ and admits the kernel of $d_A^*:H^r(M,T^*M\otimes\operatorname{ad}\mathbf P)\to H^{r-1}(M,\operatorname{ad}\mathbf P)$ as a complement.
The sum is orthogonal for the $L^2(M,T^*M\otimes\operatorname{ad}\mathbf P)$ inner product, and the orbit of $A$
under $\mathcal G^{r+1}_{\bullet}(A_0)$ has an equivariant tube and a
slice. The quotient of the component by
$\mathcal G_\bullet^{r+1}(A_0)$ has a regular quasistratification by
symmetry types. These types use the full stabilizers $S(A)$ of
Proposition~\ref{prop:estabilizadores-completos-fijos-gauge} and the
conjugacy classes
\[
 (S):=\{\mathbf u^{-1}S\mathbf u
                  \mid\mathbf u\in\mathcal G_\bullet^{r+1}(A_0)\}.
\]
The strata are the nonempty sets

\[
 \mathcal R_{(S)}:=\{[A]\in\mathcal A^r_{A_0}(\mathbf P)/
                       \mathcal G_\bullet^{r+1}(A_0)\mid S(A)\in(S)\}.
\]
Each stratum carries a manifold structure and is locally closed.
Countability of the family is not required. Regularity means that any two
strata $\Omega,\Omega'$ satisfy

\[
 \overline{\Omega}\cap\Omega'\ne\varnothing
 \quad\Longrightarrow\quad
 \begin{cases}
 \Omega'\subseteq\overline{\Omega},\\
 \overline{\Omega'}\cap(\Omega\cup\Omega')=\Omega'.
 \end{cases}
\]
Closures are relative to the quotient. The type is well defined because
$S(A\triangleleft\mathbf u)=\mathbf u^{-1}S(A)\mathbf u$.
The family of sets $\mathcal C_{(S)}$ has the same regular
quasistratification properties in the space of connections.
\end{theorem}

\begin{proof}
Proposition~\ref{prop:persistencia-espectral-componente-gauge}
extends the gap from $A_0$ to the entire component. The Hodge
decomposition and bounds for its projections follow from
Proposition~\ref{prop:coulomb-sobolev-brecha-uniforme};
Corollary~\ref{cor:carta-coulomb-sobolev-brecha} supplies Coulomb
charts, and Proposition~\ref{prop:tubo-calibre-brecha-conexa}
constructs equivariant tubes. Corollary~\ref{cor:cociente-hilbert-brecha-gauge}
gives the quotient structure, and
Proposition~\ref{prop:estratos-hilbert-estabilizadores-gauge} gives
the submanifolds of fixed type, their closures, and the submersions
onto strata. It remains to verify the two incidence conditions.

Let $\mathcal J$ be the set of stabilizer types realized in the
component. For $(S),(T)\in\mathcal J$, define
\[
 (S)\preccurlyeq(T)
 \quad\Longleftrightarrow\quad
 \text{there is }\mathbf u\in\mathcal G_\bullet^{r+1}(A_0)
                \text{ such that }\mathbf u^{-1}S\mathbf u\subseteq T.
\]
This relation is independent of representatives: replacing them by
conjugates replaces the conjugating elements by the corresponding
products. It is reflexive and transitive.
It is also antisymmetric. To prove this, suppose that inclusions hold
in both directions, and evaluate the stabilizers at the same
$p\in\mathbf P$. We obtain inclusions, up to conjugation in $G$,
between the compact Lie subgroups
$H_S=\operatorname{ev}_p(S)$ and $H_T=\operatorname{ev}_p(T)$.
The opposite inclusions first force $\dim H_S=\dim H_T$.
In each inclusion, the identity components agree, since they have
the same Lie algebra. The induced inclusions between the finite
sets of components then give the same number of components.
Thus the first inclusion of compact groups is an equality: it
contains the identity component and all components of the larger
group. Evaluation is injective on stabilizers, so the original
inclusion $\mathbf u^{-1}S\mathbf u\subseteq T$ is also an equality.
Hence $(S)=(T)$.

For $B\in\mathcal C_{(T)}$, the definition of saturation gives
\begin{equation}
\label{eq:orden-tipos-pertenencia-saturacion-gauge}
 B\in\mathcal D_S
 \quad\Longleftrightarrow\quad (S)\preccurlyeq(T).
\end{equation}
Thus, if $\mathcal D_S$ contains one point of
$\mathcal C_{(T)}$, it contains all of $\mathcal C_{(T)}$.
Proposition~\ref{prop:clausura-saturaciones-estabilizadores-gauge}
implies
\[
 \overline{\mathcal C_{(S)}}\cap\mathcal C_{(T)}\ne\varnothing
 \quad\Longrightarrow\quad
 \mathcal C_{(T)}\subseteq\overline{\mathcal C_{(S)}}.
\]
If $(S)\ne(T)$ and this intersection occurs, antisymmetry prevents
$\overline{\mathcal C_{(T)}}$ from meeting $\mathcal C_{(S)}$.
Consequently,
\[
 \overline{\mathcal C_{(T)}}\cap
           (\mathcal C_{(S)}\cup\mathcal C_{(T)})=\mathcal C_{(T)}.
\]
For $(S)=(T)$ the same equality is immediate. The sets of fixed type
are nonempty by the choice of $\mathcal J$, disjoint by definition of
type, and cover the component. This proves the regular
quasistratification of the space of connections.

On the quotient, identities
\eqref{eq:clausura-estrato-cociente-gauge} and
\eqref{eq:orden-tipos-pertenencia-saturacion-gauge} give precisely
\[
 \overline{\mathcal R_{(S)}}\cap\mathcal R_{(T)}\ne\varnothing
 \quad\Longleftrightarrow\quad (S)\preccurlyeq(T).
\]
The forward implication follows by taking a representative in
$\mathcal D_S$ of the intersection point; conversely, every
representative of type $(T)$ belongs to $\mathcal D_S$.
Inclusion of the entire stratum and the second incidence condition
follow as above from transitivity and antisymmetry.
The quotient partition is exhaustive because each orbit has a
well-defined type. All assertions are proved.
\end{proof}

The gauge tube corresponds to Theorem 4.17, pp.~283--284, of
\cite{EichhornHeberGauge}. The properties of fixed sets, exact-type
sets, and their saturations are stated in results 5.1--5.8 of that
reference. The quasistratifications of the space of connections and
of the quotient are its Theorems 5.9 and 5.13, pp.~285--286; see also
Theorem 1.3 in Chapter IV, p.~445, of \cite{Eichhorn2007}.
The types used here are defined by the full stabilizers $S(A)$;
isotropy in the selected group is $C_\bullet$, as stated in
\eqref{eq:estabilizador-grupo-seleccionado-calibre}.

\pagebreak
There is a parallel result for the action of diffeomorphisms on metrics.
For a metric $\mathbf{g}$, denote by
$\Delta_{\mathbf{g}}^{(0)}=d_{\mathbf{g}}^*d$ the nonnegative self-adjoint realization of the Beltrami
Laplacian on $L^2(M,d\lambda_{\mathbf{g}})$.

\begin{theorem}[Eichhorn: slice for metrics of negative curvature]
\label{teo:eichhorn-slice-metricas}
Let $(M,\mathbf{g}_0)$ be a smooth noncompact Riemannian manifold of dimension $m$
without boundary, of bounded geometry of all orders, with $m\geq2$, and let
$r\in\mathbb N$ with $r>\frac{m}{2}+2$. Suppose there is $c>0$ such that
$K_{\mathbf{g}_0}\leq-c$. Let $\mathbf{g}$ be a smooth metric in the Sobolev
component of order $r$ of $\mathbf{g}_0$, with bounded geometry of all orders and
$K_{\mathbf g}\leq-c_{\mathbf g}$ for some constant $c_{\mathbf g}>0$.
Require, in addition,
the positive equality of gaps
\[
\inf\sigma_{\mathrm e}(\Delta_{\mathbf{g}_0}^{(0)})
=\inf\sigma_{\mathrm e}(\Delta_{\mathbf{g}}^{(0)})>0.
\]
Consider the pullback action of the identity component of the
diffeomorphism group of order $r+1$ on the metrics of uniformly negative
sectional curvature in the Sobolev component of $\mathbf g_0$.
The action on that entire region is proper, and the orbit of
$\mathbf g$ admits a topological slice with an invariant open tube.
The action is continuous at the indicated Sobolev levels; differentiating
it requires the level changes associated with composition and pullback.
\end{theorem}

This assertion is contained in Theorem~5.8 of
\cite[p.~115]{EichhornSlice}; see also Theorem~4.17 in Chapter IV,
p.~501, of \cite{Eichhorn2007}. The harmonic construction in
Theorem~\ref{teo:rebanada-armonica-metrica-referencia} and the change of
reference in Corollary~\ref{cor:rebanada-metrica-cambio-referencia-uniforme}
provide the tube of that orbit. Properness on the entire region is
proved in Theorem~\ref{teo:accion-propia-region-negativa-completa}.

\begin{lemma}[Isometries with bounded displacement]
\label{lem:rigidez-isometria-desplazamiento-acotado}
Let $(N,\mathbf h)$ be a smooth complete connected simply connected Riemannian
manifold of dimension $m$ without boundary, with nonpositive sectional curvature and
strictly negative Ricci tensor. If $\phi:N\to N$ is an isometry and
there is $L<\infty$ such that
$d_{\mathbf h}(x,\phi(x))\leq L$ for every $x\in N$, then
$\phi=\operatorname{id}_N$.
\end{lemma}

\begin{proof}
Cartan--Hadamard, Theorem~\ref{teo:cartan-hadamard}, provides a unique
geodesic between each pair of points and smooth dependence of that
geodesic on its endpoints. Define
\[
 u(x):=\tfrac12d_{\mathbf h}(x,\phi(x))^2,
 \qquad x\in N.
\]
This function is smooth, including where $\phi(x)=x$, since squared
distance is the squared norm of the exponential logarithm.
Let $c:\mathbb R\to N$ be a geodesic. For $s\in\mathbb R$ and
$0\leq t\leq1$, let $\alpha(s,t)$ be the affinely parametrized
geodesic from $c(s)$ to $\phi(c(s))$.
Set $T=\partial_t\alpha$ and $J=\partial_s\alpha$.
The second variation formula in
Proposition~\ref{prop:segunda-variacion-energia-forma-indice} gives
\begin{equation}
\label{eq:convexidad-desplazamiento-isometria}
 \frac{d^2}{ds^2}u(c(s))
 =\int_0^1\left(|D_tJ|_{\mathbf h}^2
          -\langle R^{\mathbf h}(J,T)T,J\rangle_{\mathbf h}\right)dt
 \geq0.
\end{equation}
The boundary term vanishes because the two endpoint curves,
$c$ and $\phi\circ c$, are geodesics. Each integrand is nonnegative
by the curvature hypothesis.
The function $u\circ c$ is convex on all of $\mathbb R$ and bounded
by $L^2/2$, so it is constant. Indeed, a positive derivative would
force unbounded growth toward $+\infty$, and a negative derivative
would force it toward $-\infty$. Since every pair of points of $N$
can be joined by a geodesic, $u$ is constant.

Suppose that this constant is positive. Fix $x\in N$ and an
orthonormal basis $e_1,\ldots,e_m$ of $T_xN$. For each
$i=1,\ldots,m$, use the geodesic $c_i(s)=\exp_x(se_i)$ in
\eqref{eq:convexidad-desplazamiento-isometria}, and write $J_i$ for
the variational field along the geodesic from $x$ to $\phi(x)$.
Its velocity $T$ is nonzero and $J_i(0)=e_i$.
Since $u$ is constant, the integral of each nonnegative summand is
zero. Continuity then gives, at $t=0$,
$\langle R^{\mathbf h}(e_i,T(0))T(0),e_i\rangle_{\mathbf h}=0$.
Summing yields
\[
 \operatorname{Ric}_{\mathbf h}(T(0),T(0))
  =\sum_{i=1}^{m}
      \langle R^{\mathbf h}(e_i,T(0))T(0),e_i\rangle_{\mathbf h}=0,
\]
contradicting strict negativity of Ricci.
Thus $u=0$ and $\phi$ is the identity.
\end{proof}

\begin{proposition}[Absence of isotropy in the Sobolev component]
\label{prop:accion-difeomorfismos-negativa-libre}
Let $(M,\mathbf g_0)$ be a smooth noncompact Riemannian
manifold of dimension $m$ without boundary, with bounded geometry of order $k$.
Let $k,r\in\mathbb N$ satisfy $k>r+1>m/2+3$, and let $\mathbf g$
be a smooth metric uniformly equivalent to $\mathbf g_0$, with
$K_{\mathbf g}\leq0$ and $\operatorname{Ric}_{\mathbf g}<0$.
If $f\in\mathcal D^{r+1}_0(M,\mathbf g_0)$ satisfies
$f^*\mathbf g=\mathbf g$, then $f=\operatorname{id}_M$.
\end{proposition}

\begin{proof}
Uniform equivalence makes $\mathbf g$ complete.
Fix $\Lambda\geq1$ such that
$\Lambda^{-2}\mathbf g_0\leq\mathbf g\leq\Lambda^2\mathbf g_0$.
The exponential charts of
Theorem~\ref{teo:eichhorn-variedades-mapeos-abiertas} and openness of the
diffeomorphism group in Theorem~\ref{teo:eichhorn-difeomorfismos-abiertos}
provide a symmetric identity neighborhood $\mathcal W$ whose
elements have the form
\[
 f_X(x)=\exp_x^{\mathbf g_0}X(x),\qquad X\in H^{r+1}(M,TM).
\]
The local construction also applies to a nonorientable manifold,
since it uses exponential charts and control of $df$.
We may shrink the neighborhood until $f_{tX}$ is a diffeomorphism for every
$0\leq t\leq1$. The identity component is generated by
$\mathcal W$: the subgroup generated is open, its cosets are open,
and it is therefore also closed. Since it contains the identity
and is contained in its component, the two coincide.
Thus there are an integer $q\geq1$ and
$X_1,\ldots,X_q\in H^{r+1}(M,TM)$ such that
\[
 f=f_{X_q}\circ\cdots\circ f_{X_1}.
\]
Each $X_j$ is uniformly bounded by Sobolev embedding.

Work on one connected component $M_\alpha$ of $M$, preserved by all these
diffeomorphisms, and let $\pi:\widetilde M_\alpha\to M_\alpha$ be its universal
cover. The homotopy $f_{tX_j}$ lifts from the identity and determines
a diffeomorphism $\widetilde f_j$.
The lifted path of each point has length, for
$\widetilde{\mathbf g}=\pi^*\mathbf g$, at most
$\Lambda\|X_j\|_{C_b^0(M,TM;\mathbf g_0)}$.
Thus, for $\widetilde f=\widetilde f_q\circ\cdots\circ\widetilde f_1$,
the triangle inequality at the successive points gives
\begin{equation}
\label{eq:desplazamiento-levantado-difeomorfismo-sobolev}
 d_{\widetilde{\mathbf g}}(z,\widetilde f(z))
 \leq\Lambda\sum_{j=1}^{q}
            \|X_j\|_{C_b^0(M,TM;\mathbf g_0)},
 \qquad z\in\widetilde M_\alpha.
\end{equation}

The map $\widetilde f$ is an isometry of
$\widetilde{\mathbf g}$ because $f^*\mathbf g=\mathbf g$.
Sobolev regularity initially gives class $C^3$. The local identity
between an isometry and exponential maps,
$f(\exp_x^{\mathbf g}v)=\exp_{f(x)}^{\mathbf g}(df_xv)$,
also shows that it is smooth. The lifted metric is complete, has
nonpositive sectional curvature, and has strictly negative Ricci
tensor. Lemma~\ref{lem:rigidez-isometria-desplazamiento-acotado} and
\eqref{eq:desplazamiento-levantado-difeomorfismo-sobolev} give
$\widetilde f=\operatorname{id}_{\widetilde M_\alpha}$.
Descending and repeating on each component yields
$f=\operatorname{id}_M$.
\end{proof}

This absence of isotropy is the Bochner-type conclusion of
Theorem 5.7, p.~114, in \cite{EichhornSlice}. Its proof uses the
nonpositive curvature geometry and membership in the identity
component; it does not require the spectral condition used for
closedness of orbits and the tube in the slice theorem.

\pagebreak

Fix a smooth Riemannian manifold $(M,\mathbf g)$ (of dimension $m$) without
boundary and of bounded geometry of all orders, and an integer
$r>m/2+2$. We write
\[
 \mathscr M^r:=\operatorname{Met}^{2,r}_{\mathbf g}(M),\qquad
 \mathscr D^{r+1}:=\mathcal D^{r+1}_0(M,\mathbf g).
\]
For each tensor bundle $\mathbf F\to M$, the $H^j(M,\mathbf F;\mathbf g)$ norms are based on $\mathbf g$.
The connection of a metric $\mathbf h\in\mathscr M^r$ and its
derivatives are interpreted weakly whenever the order requires it.

\begin{proposition}[Pullback and an invariant Sobolev distance]
\label{prop:distancia-sobolev-invariante-metricas}
Pullback defines a continuous action of $\mathscr D^{r+1}$ on
$\mathscr M^r$. For $\mathbf h\in\mathscr M^r$ and
$\mathbf a\in H^r(M,S^2T^*M)$, set
\begin{equation}
\label{eq:norma-fuerte-sobolev-metricas}
 |\mathbf a|_{r,\mathbf h}
 :=\|\mathbf a\|_{H^r(M,S^2T^*M;\mathbf h)}.
\end{equation}
This family of norms defines a length distance $d_r$ inducing
the $H^r(M,S^2T^*M)$ topology and satisfying
\[
 d_r(f^*\mathbf h,f^*\mathbf k)=d_r(\mathbf h,\mathbf k),
 \qquad f\in\mathscr D^{r+1},\quad \mathbf h,\mathbf k\in\mathscr M^r.
\]
\end{proposition}

\begin{proof}
In coordinates, for $1\leq i,j\leq m$,
\[
 (f^*\mathbf a)_{ij}
 =\sum_{b,c=1}^{m}(a_{bc}\circ f)\,\partial_i f^b\,\partial_j f^c.
\]
The composition and multiplication results used in
Theorem~\ref{teo:eichhorn-difeomorfismos-abiertos} show that
$(f,\mathbf a)\mapsto f^*\mathbf a$ is continuous with values in $H^r(M,S^2T^*M)$.
Indeed, composition preserves $H^r(M,S^2T^*M)$, the factors $df$ have local
$H^r(U,\mathbb R^{m\times m})$ regularity. In the charts, their reference parts have uniformly
bounded derivatives through order $r$, and the differences belong
to $H^r(U,\mathbb R^{m\times m})$; both act as multipliers.
The constants are controlled in a neighborhood
of each $f$ by its bi-Lipschitz bounds and Sobolev coordinates.
The preceding formula and change of variables give the same
assertion for the global norms.

The reference term is treated separately. In an identity chart,
$f_X(x)=\exp_x^{\mathbf g}X(x)$, the expression
$f_X^*\mathbf g-\mathbf g$ is a smooth function of $X$ and $\nabla X$
that vanishes when both are zero. The exponential jet bounds and
the superposition results already proved give a smooth map from
a ball in $H^{r+1}(M,TM)$ to $H^r(M,S^2T^*M)$.
For an arbitrary element of the identity component, use a finite
factorization into elements of that chart. The identity
\[
 f^*\mathbf h-\mathbf g
  =(f^*\mathbf g-\mathbf g)+f^*(\mathbf h-\mathbf g)
\]
proves membership in the affine component and joint continuity.
Uniform positivity is preserved because $f$ and $f^{-1}$ have
uniformly bounded differentials.

Expansion \eqref{eq:expansion-cambio-conexion-tensor-metricas}
and the Sobolev products show that \eqref{eq:norma-fuerte-sobolev-metricas}
is finite and depends continuously on $(\mathbf h,\mathbf a)$.
Apply the expansion first to smooth metrics; approximation in $H^r(M,S^2T^*M)$
with a common positive lower bound gives the assertion for
$\mathbf h\in\mathscr M^r$.
More precisely, for every $\mathbf h_*$ there are $\rho>0$ and
$0<a\leq A<\infty$ such that
\begin{equation}
\label{eq:comparacion-local-norma-fuerte-metricas}
 a\|\mathbf a\|_{H^r(M,S^2T^*M)}
 \leq |\mathbf a|_{r,\mathbf h}
 \leq A\|\mathbf a\|_{H^r(M,S^2T^*M)},\qquad
 \|\mathbf h-\mathbf h_*\|_{H^r(M,S^2T^*M)}<\rho.
\end{equation}
Here $\rho$ is decreased to retain a positive lower bound.
The Koszul formula controls
$\nabla^{\mathbf h}-\nabla^{\mathbf g}$ in $H^{r-1}(M,T^*M\otimes\operatorname{End}TM)$, uniformly
on this ball; the cited expansion and its reverse provide the two
constants. No high-order curvature bounds are required for the
variable metrics in this ball.

For a path $\gamma:[0,1]\to\mathscr M^r$ that is piecewise $C^1$, define
\[
 L_r(\gamma):=\int_0^1|\dot\gamma(t)|_{r,\gamma(t)}\,dt,
 \qquad
 d_r(\mathbf h,\mathbf k):=
 \inf_{\substack{\gamma(0)=\mathbf h,\ \gamma(1)=\mathbf k\\
                 \gamma\text{ is piecewise }C^1}}
                 L_r(\gamma).
\]
Convexity of $\mathscr M^r$ ensures that the infimum is finite.
If $\mathbf h,\mathbf k\in B_{H^r(M,S^2T^*M)}(\mathbf h_*,\rho/4)$,
the segment gives $d_r(\mathbf h,\mathbf k)\leq
A\|\mathbf h-\mathbf k\|_{H^r(M,S^2T^*M)}$.
A path remaining in the ball of radius $\rho$ has length at least
$a\|\mathbf h-\mathbf k\|_{H^r(M,S^2T^*M)}$.
If it leaves the ball, its portion up to the first exit has length
at least $3a\rho/4$, whereas
$\|\mathbf h-\mathbf k\|_{H^r(M,S^2T^*M)}<\rho/2$.
Consequently,
\[
 a\|\mathbf h-\mathbf k\|_{H^r(M,S^2T^*M)}
 \leq d_r(\mathbf h,\mathbf k)
 \leq A\|\mathbf h-\mathbf k\|_{H^r(M,S^2T^*M)}
\]
on the smaller ball. To reach a point outside that ball from
$\mathbf h_*$, every path has length at least $a\rho/4$.
This proves separation of points and local equivalence of topologies.
Symmetry and the triangle inequality
follow by reversing and concatenating paths.

Naturality of the connection and change of variables give
\[
 |f^*\mathbf a|_{r,f^*\mathbf h}=|\mathbf a|_{r,\mathbf h}.
\]
For Sobolev objects, the identity follows by approximation in the
norms already controlled. Pullback by a fixed $f$ is affine and
continuous in metric charts, and takes piecewise $C^1$ paths to
paths of the same class and length. Applying it first to $f$ and
then to $f^{-1}$ proves invariance of $d_r$.
\end{proof}

\begin{lemma}[Coercivity and operators of the metric action]
\label{lem:coercividad-operadores-rebanada-metrica}
Suppose that $\operatorname{Ric}_{\mathbf g}\leq-\kappa\mathbf g$
for some $\kappa>0$. Define
\[
 \begin{aligned}
 D_{\mathbf g}X&:=\mathcal L_X\mathbf g,\\
 B_{\mathbf g}\mathbf a&:=
 \bigl(\operatorname{div}_{\mathbf g}\mathbf a
           -\tfrac12 d\operatorname{tr}_{\mathbf g}\mathbf a\bigr)^\sharp,\\
 J_{\mathbf g}X&:=\nabla^*\nabla X-
                         \operatorname{Ric}_{\mathbf g}^{\sharp}X.
 \end{aligned}
\]
The divergence of a symmetric tensor has components
$(\operatorname{div}_{\mathbf g}\mathbf a)_j
=\displaystyle\sum_{i=1}^{m}(\nabla_{e_i}\mathbf a)(e_i,e_j)$
in an orthonormal frame. Then $J_{\mathbf g}$ and
$D_{\mathbf g}^*D_{\mathbf g}$ are isomorphisms from $H^{r+1}(M,TM)$
onto $H^{r-1}(M,TM)$, and there are constants $C_J,C_D>0$ such that
\begin{equation}
\label{eq:coercividad-sobolev-accion-metrica}
 \|X\|_{H^{r+1}(M,TM)}\leq C_J\|J_{\mathbf g}X\|_{H^{r-1}(M,TM)},
 \qquad
 \|X\|_{H^{r+1}(M,TM)}\leq C_D\|D_{\mathbf g}X\|_{H^r(M,S^2T^*M)}.
\end{equation}
Moreover, $-B_{\mathbf g}D_{\mathbf g}=J_{\mathbf g}$ and
\begin{equation}
\label{eq:descomposicion-tangente-metricas}
 H^r(M,S^2T^*M)=D_{\mathbf g}H^{r+1}(M,TM)
                   \oplus\ker D_{\mathbf g}^*.
\end{equation}
The last sum is orthogonal in $L^2(M,S^2T^*M)$ and has continuous projections
in $H^r(M,S^2T^*M)$.
\end{lemma}

\begin{proof}
For smooth compactly supported $X$, integration by parts and
commutation of two covariant derivatives give
\begin{equation}
\label{eq:identidad-korn-ricci-metricas}
 \|D_{\mathbf g}X\|_{L^2(M,S^2T^*M)}^2
 =2\int_M\bigl(|\nabla X|_{\mathbf g}^2
       +(\operatorname{div}_{\mathbf g}X)^2
       -\operatorname{Ric}_{\mathbf g}(X,X)\bigr)
          \,d\lambda_{\mathbf g}.
\end{equation}
Indeed, expanding the square of
$(D_{\mathbf g}X)_{ij}=\nabla_iX_j+\nabla_jX_i$, the cross term
is twice the integral of
$\displaystyle\sum_{i,j=1}^{m}\nabla_iX_j\nabla_jX_i$.
Integrating first in $i$, commuting the derivatives, and then
integrating in $j$ transforms this integral into
\[
 \int_M\bigl((\operatorname{div}_{\mathbf g}X)^2
                -\operatorname{Ric}_{\mathbf g}(X,X)\bigr)
              \,d\lambda_{\mathbf g}.
\]
Density of smooth compactly supported fields in $H^1(M,TM)$ extends
\eqref{eq:identidad-korn-ricci-metricas} to that space. In particular,
\[
 \|D_{\mathbf g}X\|_{L^2(M,S^2T^*M)}^2
 \geq2\|\nabla X\|_{L^2(M,T^*M\otimes TM)}^2+2\kappa\|X\|_{L^2(M,TM)}^2.
\]
The same calculation identifies
\[
 D_{\mathbf g}^*\mathbf a=-2(\operatorname{div}_{\mathbf g}\mathbf a)^\sharp,
 \quad
 D_{\mathbf g}^*D_{\mathbf g}
   =2J_{\mathbf g}-2\operatorname{grad}_{\mathbf g}\operatorname{div}_{\mathbf g},
 \quad -B_{\mathbf g}D_{\mathbf g}=J_{\mathbf g}.
\]
The principal symbol of the second operator, applied to $v\in T_xM$,
is
\[
 \sigma_2(D_{\mathbf g}^*D_{\mathbf g})(x,\xi)v
    =2\bigl(|\xi|_{\mathbf g}^2v+\xi(v)\xi^\sharp\bigr).
\]
Its inner product with $v$ is at least
$2|\xi|_{\mathbf g}^2|v|_{\mathbf g}^2$; the symbol of
$J_{\mathbf g}$ is $|\xi|_{\mathbf g}^2I$.

The bilinear forms
\[
 \begin{aligned}
 q_J(X,Y)&:=\int_M
   \bigl(\langle\nabla X,\nabla Y\rangle_{\mathbf g}
                 -\operatorname{Ric}_{\mathbf g}(X,Y)\bigr)
                 \,d\lambda_{\mathbf g},\\
 q_D(X,Y)&:=\int_M
         \langle D_{\mathbf g}X,D_{\mathbf g}Y\rangle_{\mathbf g}
                 \,d\lambda_{\mathbf g}
 \end{aligned}
\]
are continuous and coercive on $H^1(M,TM)$, with lower constants
$\min\{1,\kappa\}$ and $2\min\{1,\kappa\}$.
The Riesz representation theorem, applied to each of these
equivalent Hilbert norms, gives, for every $F\in L^2(M,TM)$, a unique
weak solution $X\in H^1(M,TM)$ of $J_{\mathbf g}X=F$ or
$D_{\mathbf g}^*D_{\mathbf g}X=F$, respectively.
Taking $Y=X$ gives
$\|X\|_{L^2(M,TM)}\leq\kappa^{-1}\|F\|_{L^2(M,TM)}$ in the first case and
$\|X\|_{L^2(M,TM)}\leq(2\kappa)^{-1}\|F\|_{L^2(M,TM)}$ in the second.

The localization used in
Proposition~\ref{prop:coulomb-sobolev-brecha-uniforme} applies to
these two smooth uniformly elliptic operators. In synchronous
frames, their coefficients and the required derivatives are
uniformly bounded. The second-order estimate is obtained by
freezing the preceding symbol; its eigenvalues lie between
$2|\xi|^2$ and $4|\xi|^2$ for $D_{\mathbf g}^*D_{\mathbf g}$.
After differentiating through order $r-1$, the commutators contain
derivatives of $X$ of order at most $r$. The cutoffs and bounded
multiplicity of the cover give, for $P=J_{\mathbf g}$ or
$P=D_{\mathbf g}^*D_{\mathbf g}$,
\[
 \|X\|_{H^{r+1}(M,TM)}
 \leq C_P\bigl(\|PX\|_{H^{r-1}(M,TM)}+\|X\|_{L^2(M,TM)}\bigr).
\]
Here the intermediate orders are absorbed using
$\|X\|_{H^r(M,TM)}\leq\varepsilon\|X\|_{H^{r+1}(M,TM)}+
C_\varepsilon\|X\|_{L^2(M,TM)}$, with sufficiently small $\varepsilon>0$.
Local elliptic regularity and these global estimates allow their
application to weak solutions. This proves both isomorphisms.
Only geometric bounds through order $r+1$ enter this step.

Write $G_D=(D_{\mathbf g}^*D_{\mathbf g})^{-1}$.
The bound for $G_D$ and continuity of $D_{\mathbf g}^*$ give the
second inequality in \eqref{eq:coercividad-sobolev-accion-metrica}.
Thus $D_{\mathbf g}$ is injective with closed range.
For every $\mathbf a\in H^r(M,S^2T^*M)$, the decomposition is
\[
 \mathbf a=D_{\mathbf g}G_DD_{\mathbf g}^*\mathbf a
       +\bigl(\mathbf a-D_{\mathbf g}G_DD_{\mathbf g}^*\mathbf a\bigr).
\]
The second term belongs to the kernel of $D_{\mathbf g}^*$.
Both operators are continuous in $H^r(M,S^2T^*M)$, and the definition of
the adjoint proves orthogonality and uniqueness.
\end{proof}

For a metric $\mathbf h\in\mathscr M^r$ and a map
$f\in\mathscr D^{r+1}$, its \emph{tension field} with target metric
$\mathbf g$ is
\[
 \tau_{\mathbf h,\mathbf g}(f)
       :=\operatorname{tr}_{\mathbf h}\nabla df.
\]
The connection in $\nabla df$ uses $\mathbf h$ on the domain and
$\mathbf g$ on the target. We call $f$ \emph{harmonic} when this
field vanishes. The definition is local and does not require
finite total energy of the map.

\begin{proposition}[Local harmonic gauge for metrics]
\label{prop:calibre-armonico-local-metricas}
If $\operatorname{Ric}_{\mathbf g}\leq-\kappa\mathbf g$ with
$\kappa>0$, there are a neighborhood $\mathcal O$ of $\mathbf g$
in $\mathscr M^r$, a number $b>0$, and a smooth map
\[
 X:\mathcal O\longrightarrow B_{H^{r+1}(M,TM)}(0,b),
 \qquad X(\mathbf g)=0,
\]
such that $F(\mathbf h):=f_{X(\mathbf h)}$ is a harmonic
diffeomorphism from $(M,\mathbf h)$ to $(M,\mathbf g)$.
It is the unique solution of this form with $X$ in that ball.
Moreover, near $\mathbf g$, the set
\begin{equation}
\label{eq:conjunto-calibre-armonico-metricas}
 \mathcal S_0:=\{\mathbf h\in\mathcal O\mid
                         \tau_{\mathbf h,\mathbf g}(\operatorname{id})=0\}
\end{equation}
is a smooth Hilbert submanifold, with
$T_{\mathbf g}\mathcal S_0=\ker B_{\mathbf g}$, and
\[
 H^r(M,S^2T^*M)=D_{\mathbf g}H^{r+1}(M,TM)\oplus\ker B_{\mathbf g}.
\]
\end{proposition}

\begin{proof}
Take $b$ small enough that all $f_X$ with
$\|X\|_{H^{r+1}(M,TM)}<b$ belong to the exponential chart of
$\mathscr D^{r+1}$. Denote by $P_X(x)$ parallel transport from
$x$ to $f_X(x)$ along
$t\mapsto\exp_x^{\mathbf g}(tX(x))$, $0\leq t\leq1$.
Consider
\[
 \mathcal T(\mathbf h,X)
    :=P_X^{-1}\tau_{\mathbf h,\mathbf g}(f_X)
       \in H^{r-1}(M,TM).
\]
We specify its regularity. In charts, for $1\leq a\leq m$,
\begin{equation}
\label{eq:tension-coordenadas-indices}
 \tau_{\mathbf h,\mathbf g}(f)^a
 =\sum_{i,j=1}^{m}h^{ij}\left(
       \partial_i\partial_j f^a
       -\sum_{k=1}^{m}\Gamma(\mathbf h)^k_{ij}\partial_k f^a
       +\sum_{b,c=1}^{m}\Gamma(\mathbf g)^a_{bc}(f)
                     \partial_i f^b\partial_j f^c\right).
\end{equation}
In uniform normal charts, the exponential and parallel transport
have the bounds of Lemma~\ref{lem:jets-mixtos-exponencial-orden-finito}.
Decrease $b$ so that all points
and geodesics used remain in larger charts of a common radius.
After substituting $f_X$, the expression depends smoothly on
$X,\nabla X,\mathbf h-\mathbf g$, and affinely on the second
derivatives of $X$ and first derivatives of $\mathbf h-\mathbf g$.
The coefficients and their derivatives are uniformly bounded on
small balls of these variables. Inversion of $\mathbf h$ retains
these properties when a positive lower bound is imposed.
Since $r-1>m/2$, Sobolev products and superposition show that
$\mathcal T$ is smooth. In particular, for every integer
$\ell\geq1$, there is $C_\ell>0$, on a fixed neighborhood of
$(\mathbf g,0)$, such that
\[
 \left\|D^\ell\mathcal T(\mathbf h,X)
           [(\mathbf a_1,Y_1),\ldots,(\mathbf a_\ell,Y_\ell)]
       \right\|_{H^{r-1}(M,TM)}
 \leq C_\ell\prod_{i=1}^{\ell}
              \bigl(\|\mathbf a_i\|_{H^r(M,S^2T^*M)}+\|Y_i\|_{H^{r+1}(M,TM)}\bigr).
\]
These bounds follow by differentiating the jet formula: each term
contains factors in $H^{r-1}(M,\mathbf F)$ or a higher order and uniform
coefficients. The zeroth-order term is controlled by subtracting
$\mathcal T(\mathbf g,0)=0$ and using the integral formula for
differences. An integrable factor therefore always remains.
The same integral formula applied to the derivatives, as in the
superposition proof, shows that they are Fr\'echet differentials.

The variation of the Levi--Civita connection, obtained from
\eqref{eq:diferencia-levi-civita-metricas}, and commutation of
derivatives in the variation of $f$ give
\begin{equation}
\label{eq:linealizacion-tension-metricas}
 D_{\mathbf h}\mathcal T(\mathbf g,0)[\mathbf a]=-B_{\mathbf g}\mathbf a,
 \qquad D_X\mathcal T(\mathbf g,0)[Y]=-J_{\mathbf g}Y.
\end{equation}
For the second identity, the second-order term is
$\operatorname{tr}_{\mathbf g}\nabla^2Y=-\nabla^*\nabla Y$;
commutation adds $\operatorname{Ric}_{\mathbf g}^{\sharp}Y$.
Derivatives of parallel transport do not contribute, since the
tension at $(\mathbf g,\operatorname{id})$ is zero.
Lemma~\ref{lem:coercividad-operadores-rebanada-metrica} and
Theorem~\ref{teo:funcion-implicita-banach} provide $X(\mathbf h)$,
its smoothness, and uniqueness in a smaller ball.
We may take $\mathcal O$ to be a sufficiently small convex ball;
the path $t\mapsto f_{X(\mathbf g+t(\mathbf h-\mathbf g))}$ joins
the identity to $F(\mathbf h)$ in the indicated group.

Finally, $\mathbf h\mapsto\mathcal T(\mathbf h,0)$ has differential
at $\mathbf g$ equal to $-B_{\mathbf g}$, with right inverse
$D_{\mathbf g}J_{\mathbf g}^{-1}$, since
$-B_{\mathbf g}D_{\mathbf g}=J_{\mathbf g}$.
The implicit function theorem proves the assertion about
$\mathcal S_0$. Projection onto the first summand of the tangent
decomposition is $-D_{\mathbf g}J_{\mathbf g}^{-1}B_{\mathbf g}$;
its complement has range $\ker B_{\mathbf g}$.
\end{proof}

\begin{lemma}[Uniqueness of harmonic diffeomorphisms]
\label{lem:unicidad-armonica-componente-sobolev}
Suppose $m\geq2$ and $K_{\mathbf g}\leq-c$ for some $c>0$.
If $\mathbf h\in\mathscr M^r$ and $f_1,f_2\in\mathscr D^{r+1}$
are harmonic from $(M,\mathbf h)$ to $(M,\mathbf g)$, then
$f_1=f_2$.
\end{lemma}

\begin{proof}
A field $X\in H^{r+1}(M,TM)$ tends uniformly to zero outside compact
sets. The argument is the same as in the proof of
Theorem~\ref{teo:eichhorn-espacio-metricas}: if $|X|\geq\varepsilon$
along a sequence escaping every compact set, the bound for
$\nabla X$ gives disjoint balls of a common radius on which
$|X|\geq\varepsilon/2$. The lower volume bound contradicts
$X\in L^2(M,TM)$.

Every $f\in\mathscr D^{r+1}$ is a finite product
$f=f_{X_q}\circ\cdots\circ f_{X_1}$ of elements of an exponential
identity neighborhood. On each connected component of $M$, lift
their exponential isotopies to the universal covering space.
The lifts starting at the identity commute with deck transformations.
If $\widetilde f$ is the product of the lifts and $\pi$ is the
projection, the lengths of the lifted paths give
\[
 d_{\widetilde{\mathbf g}}(z,\widetilde f(z))
 \leq\sum_{j=1}^{q}|X_j(v_{j-1}(\pi(z)))|_{\mathbf g},
 \quad
 v_0=\operatorname{id},\quad
 v_{j-1}=f_{X_{j-1}}\circ\cdots\circ f_{X_1}\quad(2\leq j\leq q).
\]
The right-hand side is bounded and tends to zero outside compact
sets, since each $v_{j-1}$ is a proper homeomorphism.
Apply this construction to $f_1$ and $f_2$.
Equivariance of the lifts makes the function
\[
 u(x):=\tfrac12d_{\widetilde{\mathbf g}}
           (\widetilde f_1(z),\widetilde f_2(z))^2,
 \qquad \pi(z)=x,
\]
well defined on $M$, of class $C^2$, and tending to zero outside
compact sets. Smoothness of squared distance on the covering space
follows from Cartan--Hadamard; the maps $f_i$ are at least $C^3$
by Sobolev embedding.

If $u$ has a positive value, it attains a positive maximum at some
$x\in M$. Choose an $\mathbf h$-orthonormal basis
$e_1,\ldots,e_m$ of $T_xM$, and vary the endpoints of the geodesic
between $\widetilde f_1(z)$ and $\widetilde f_2(z)$ along the curves
$f_\nu(\exp_x^{\mathbf h}(se_i))$, for $\nu=1,2$.
Write $T$ for the velocity of this geodesic, parametrized on
$[0,1]$, and $J_i$ for its variational fields.
Summing the second variation formula of
Proposition~\ref{prop:segunda-variacion-energia-forma-indice}, the boundary terms are
inner products with the tension fields of $f_1$ and $f_2$; both
vanish. Hence
\begin{equation}
\label{eq:subarmonicidad-distancia-mapeos-armonicos}
 \operatorname{tr}_{\mathbf h}\operatorname{Hess}_{\mathbf h}u(x)
 =\sum_{i=1}^{m}\int_0^1
   \bigl(|D_tJ_i|_{\widetilde{\mathbf g}}^2
       -\langle R^{\widetilde{\mathbf g}}(J_i,T)T,J_i
                  \rangle_{\widetilde{\mathbf g}}\bigr)\,dt.
\end{equation}
All summands are nonnegative. Moreover, $T(0)\ne0$, and the vectors
$J_i(0)=d\widetilde f_1(e_i)$ form a basis.
Since $m\geq2$, at least one is not collinear with $T(0)$.
Strictly negative curvature makes the corresponding integrand
positive at $t=0$ and, by continuity, on an interval of positive
length. The right-hand side is then positive, contradicting
$\operatorname{tr}_{\mathbf h}\operatorname{Hess}_{\mathbf h}u(x)\leq0$
at a maximum. Thus $u=0$ and $f_1=f_2$.
The argument also applies to compact components, where a maximum
exists directly.
\end{proof}

\begin{theorem}[Slice at the reference metric]
\label{teo:rebanada-armonica-metrica-referencia}
Let $(M,\mathbf g)$ (of dimension $m$) be smooth, and without boundary,
of bounded geometry of all orders, with $m\geq2$ and
$K_{\mathbf g}\leq-c$ for some $c>0$.
Let $r\in\mathbb N$ satisfy $r>m/2+2$.
There are a smooth Hilbert submanifold
$\mathcal S\subseteq\mathscr M^r$ containing $\mathbf g$ and
an invariant open set $\mathcal U\subseteq\mathscr M^r$ such that
\begin{equation}
\label{eq:tubo-topologico-metricas-referencia}
 \Theta:\mathcal S\times\mathscr D^{r+1}\longrightarrow\mathcal U,
 \qquad \Theta(\mathbf s,f)=f^*\mathbf s,
\end{equation}
is a homeomorphism. The set $\mathcal U$ can be chosen to consist
of metrics with uniformly negative sectional curvature.
The action on $\mathcal U$ is proper. Every orbit contained in
$\mathcal U$ is closed in $\mathscr M^r$, and its parametrization
by $\mathscr D^{r+1}$ is a homeomorphism onto the orbit.
The orbit of $\mathbf g$ is also a smooth Hilbert submanifold,
with tangent space $D_{\mathbf g}H^{r+1}(M,TM)$ at $\mathbf g$.
\end{theorem}

\begin{proof}
We have $\operatorname{Ric}_{\mathbf g}\leq-(m-1)c\mathbf g$,
so the preceding results apply. Shrink the neighborhood
$\mathcal O$ from the harmonic-gauge proposition and define
\[
 \Pi(\mathbf h):=(F(\mathbf h)^{-1})^*\mathbf h.
\]
This map is continuous by continuity of inversion and pullback.
Naturality of the tension field under an isometry of the domain
gives $\tau_{\Pi(\mathbf h),\mathbf g}(\operatorname{id})=0$:
the identity is the composition of the harmonic map
$F(\mathbf h):(M,\mathbf h)\to(M,\mathbf g)$ with the isometry
$F(\mathbf h)^{-1}:(M,\Pi(\mathbf h))\to(M,\mathbf h)$.
Since $\Pi(\mathbf g)=\mathbf g$, take an open set
$\mathcal O_1\subseteq\mathcal O$ such that
$\Pi(\mathcal O_1)\subseteq\mathcal S_0$.

Choose an open ball $\mathcal O_0$ centered at $\mathbf g$,
contained in $\mathcal O_1$, and an open identity neighborhood
$\mathcal W$ in $\mathscr D^{r+1}$, small enough that
\[
 f^*\mathbf s\in\mathcal O_1\quad
 (\mathbf s\in\mathcal S_0\cap\mathcal O_0,\ f\in\mathcal W).
\]
Continuity at $(\mathbf g,\operatorname{id})$ ensures this choice.
Set $\mathcal S:=\mathcal S_0\cap\mathcal O_0$.
If $\mathbf h=f^*\mathbf s$ with $\mathbf s\in\mathcal S$ and
$f\in\mathcal W$, then both $f$ and $F(\mathbf h)$ are harmonic
from $(M,\mathbf h)$ to $(M,\mathbf g)$.
The uniqueness lemma gives $F(\mathbf h)=f$ and
$\Pi(\mathbf h)=\mathbf s$.
Consequently, the restriction of $\Theta$ to
$\mathcal S\times\mathcal W$ is a homeomorphism onto the open set
\[
 \mathcal V:=\{\mathbf h\in\mathcal O_1\mid
                   \Pi(\mathbf h)\in\mathcal O_0,
                   \ F(\mathbf h)\in\mathcal W\},
\]
with inverse $\mathbf h\mapsto(\Pi(\mathbf h),F(\mathbf h))$.

If $f^*\mathbf s_1=\mathbf s_2$ for
$\mathbf s_1,\mathbf s_2\in\mathcal S$ and any
$f\in\mathscr D^{r+1}$, then $f$ and the identity are harmonic
from $(M,\mathbf s_2)$ to $(M,\mathbf g)$.
The same lemma gives $f=\operatorname{id}$ and
$\mathbf s_1=\mathbf s_2$.
This proves global injectivity of \eqref{eq:tubo-topologico-metricas-referencia}.
For each $f_0\in\mathscr D^{r+1}$, the identity
\[
 \Theta(\mathbf s,f\circ f_0)=f_0^*\Theta(\mathbf s,f)
\]
transports the preceding local product to
$\mathcal S\times(\mathcal W\circ f_0)$.
These open sets cover the product; their images are the open sets
$f_0^*\mathcal V$. Thus $\Theta$ is a homeomorphism onto
\[
 \mathcal U=\bigcup_{f_0\in\mathscr D^{r+1}}f_0^*\mathcal V.
\]
The embedding $H^r(M,S^2T^*M)\hookrightarrow C_b^2(M,S^2T^*M)$ allows us to choose
$\mathcal O_0$ so that $K_{\mathbf s}\leq-c/2$ for
$\mathbf s\in\mathcal S$. Naturality of curvature preserves
this bound throughout $\mathcal U$.

In these coordinates, the map defining properness of the action
has the form
\[
 (\mathbf s,a,f)\longmapsto
          \bigl((\mathbf s,a),(\mathbf s,a\circ f)\bigr).
\]
It is a homeomorphism onto the closed subset of pairs having the
same $\mathcal S$ coordinate; its inverse uses
$f=a^{-1}\circ b$ at the pair $((\mathbf s,a),(\mathbf s,b))$.
The preimage of a compact set is therefore compact, proving
the asserted properness.

We prove closedness in the entire space $\mathscr M^r$.
Fix $\mathbf h=\Theta(\mathbf s,a)\in\mathcal U$ and a sequence
$f_n^*\mathbf h\to\mathbf k$ in $\mathscr M^r$.
Choose $\varepsilon>0$ such that
$B_{d_r}(\mathbf h,\varepsilon)\subseteq\mathcal U$, and then
$N\in\mathbb N$ with $d_r(f_N^*\mathbf h,\mathbf k)<\varepsilon$.
Invariance of $d_r$ implies
\[
 \mathbf z:=(f_N^{-1})^*\mathbf k\in\mathcal U,
 \qquad
 \mathbf z_n:=(f_N^{-1})^*f_n^*\mathbf h\longrightarrow\mathbf z.
\]
For sufficiently large $n$, the tube coordinates of $\mathbf z_n$
are $(\mathbf s,a\circ f_n\circ f_N^{-1})$.
Continuity of $\Theta^{-1}$ gives $b\in\mathscr D^{r+1}$ with
$\mathbf z=\Theta(\mathbf s,b)$ and
\[
 f_n\longrightarrow a^{-1}\circ b\circ f_N
                    \quad\text{in }\mathscr D^{r+1}.
\]
Thus $\mathbf k$ belongs to the orbit of $\mathbf h$.
Metrizability turns this sequential verification into closedness.
The same argument, or the restriction of $\Theta$ to
$\{\mathbf s\}\times\mathscr D^{r+1}$, proves the topological
assertion about the parametrization of every orbit.

For the reference orbit, $f\mapsto f^*\mathbf g$ is smooth:
at the identity it is the jet expression used in the first
proposition, and at another point it is transported by right
composition and a fixed pullback. Its differential at the identity
is $D_{\mathbf g}$, which is injective with complemented range by
\eqref{eq:descomposicion-tangente-metricas}.
Indeed, the map
$(X,\mathbf a)\mapsto f_X^*\mathbf g+\mathbf a$, with
$\mathbf a\in\ker D_{\mathbf g}^*$, has at $(0,0)$ differential
$(Y,\mathbf b)\mapsto D_{\mathbf g}Y+\mathbf b$,
which is an isomorphism.
Theorem~\ref{teo:funcion-inversa-local-banach} provides a chart in
which the local orbit corresponds to $\mathbf a=0$.
The topological parametrization of the orbit identifies these
charts with the relative topology and gives the stated submanifold.
\end{proof}

In this construction, the harmonic equation is solved as a smooth
map from $\mathscr M^r\times H^{r+1}(M,TM)$ to $H^{r-1}(M,TM)$.
Passing from those solutions to tube coordinates uses continuity
of pullback and inversion. Thus
\eqref{eq:tubo-topologico-metricas-referencia} is asserted as a
homeomorphism, while the slice and the reference orbit carry the
smooth structures described above.

\begin{lemma}[Control of a chart by the jets of the map]
\label{lem:carta-exponencial-control-jets-referencia}
Let $(M,\mathbf h)$ be a Riemannian manifold of dimension $m$, without boundary, and of bounded
geometry of all orders. Fix an integer $s>m/2+3$ and $B\geq2\sqrt m$.
There is $\rho>0$, smaller than one quarter of the injectivity
radius, with the following property. Let $f:M\to M$ have class
$C^1\cap H^s_{\mathrm{loc}}(M,M)$, with
\[
 d_{\mathbf h}(x,f(x))<\rho,\qquad |df_x|_{\mathbf h}\leq B
                         \quad(x\in M).
\]
Write $X(x)=\log_x^{\mathbf h}f(x)$ and let
$P_f(x):T_xM\to T_{f(x)}M$ be radial parallel transport. Set
\[
 K_{\mathbf h}(f):=P_f^{-1}df-I,
 \qquad Q_{\mathbf h}(f):=\nabla^{\mathbf h}df.
\]
In $Q_{\mathbf h}(f)$ we use the base connection and the induced
connection on $f^*TM$. Define
\begin{equation}
\label{eq:datos-jets-carta-referencia}
 \begin{split}
 J_{s,\mathbf h}(f):={}&
 \|d_{\mathbf h}(\,\cdot\,,f)\|_{L^2(M,\mathbf h)}
 +\|K_{\mathbf h}(f)\|_{L^2(M,\operatorname{End}TM;\mathbf h)}\\
 &+\sum_{j=0}^{s-2}
 \|\nabla^jQ_{\mathbf h}(f)\|_{L^2(M,(T^*M)^{\otimes(j+2)}
                                  \otimes f^*TM;\mathbf h)}.
 \end{split}
\end{equation}
If this number is finite, then $X\in H^s(M,TM;\mathbf h)$.
There is a polynomial $P_{s,B}$ with nonnegative coefficients and
$P_{s,B}(0)=0$ such that
\begin{equation}
\label{eq:control-logaritmo-por-jets}
 \|X\|_{H^s(M,TM;\mathbf h)}\leq P_{s,B}(J_{s,\mathbf h}(f)).
\end{equation}
Conversely, on a sufficiently small Sobolev ball, the map
$X\mapsto f=\operatorname{Exp}_{\operatorname{id}}^{\mathbf h}X$
makes the data in \eqref{eq:datos-jets-carta-referencia} continuous,
after transporting the fibers by $P_f^{-1}$, and all tend to zero
when $X\to0$ in $H^s(M,TM;\mathbf h)$.
\end{lemma}

\begin{proof}
We use the jet bounds for the exponential, logarithm and parallel
transport from Lemma~\ref{lem:jets-mixtos-exponencial-orden-finito}.
Reduce $\rho$ so that the vertical differential of the exponential
and its inverse have norm at most two. The first derivative,
decomposed horizontally and vertically, gives
\[
 P_f^{-1}df=B_X+A_X\nabla X,
 \qquad A_0=B_0=I,\qquad |B_X-I|\leq C|X|.
\]
The coefficients have uniformly bounded mixed derivatives.
Since $|X|=d_{\mathbf h}(\,\cdot\,,f)$, we obtain
\[
 |\nabla X|\leq2|K_{\mathbf h}(f)|+2C|X|,
 \qquad \|X\|_{H^1(M,TM;\mathbf h)}\leq C_1J_{s,\mathbf h}(f).
\]
Moreover, $|df|\leq B$ gives a uniform bound for $|\nabla X|$.
The norms of $X$ in this proof are over $M$ with $\mathbf h$.

Differentiating once more and solving for the second derivative gives
\begin{equation}
\label{eq:segunda-derivada-logaritmo-jets}
 \nabla^2X=A_X^{-1}P_f^{-1}Q_{\mathbf h}(f)
                         +R(x,X,\nabla X).
\end{equation}
The expression $R$ is a polynomial of degree at most two in
$\nabla X$, with smooth coefficients in $(x,X)$ and uniformly
bounded derivatives; it satisfies $R(x,0,0)=0$.
Indeed, on differentiating $P_f(B_X+A_X\nabla X)$, the only term
with two derivatives of $X$ is $P_fA_X\nabla^2X$.
The derivatives of $A_X,B_X$ are affine in $\nabla X$.
The derivative of parallel transport is given by
\eqref{eq:variacion-conexion-transporte-curvatura}, whose integrand
contains $X$ and is affine in $\nabla X$. Setting $X=0$ gives
$f=\operatorname{id}$ and $Q_{\mathbf h}(f)=0$, proving the
stated vanishing.

Let us specify the product estimates. On the balls of a uniform
cover, with fixed cutoffs, we use
\[
 \|uv\|_{H^j(\mathbb R^m)}\leq C_j\bigl(
       \|u\|_{L^\infty(\mathbb R^m)}\|v\|_{H^j(\mathbb R^m)}
       +\|v\|_{L^\infty(\mathbb R^m)}\|u\|_{H^j(\mathbb R^m)}\bigr),\qquad j\in\mathbb N.
\]
This bound follows from Leibniz, Hölder and interpolation of
intermediate derivatives between $L^\infty(\mathbb R^m)$ and $H^j(\mathbb R^m)$: a factor
with $a$ derivatives is assigned exponent $2j/a$, for
$1\leq a\leq j$. The fractions $a/(2j)$ and $(j-a)/(2j)$
sum to $1/2$. Cutoffs introduce only lower orders.
Finite multiplicity allows the local bounds to be summed.
The same rule and Taylor's formula in $X$ control the composed
coefficients, separating their value at $X=0$.

Set $q=s-2$ and $Q=Q_{\mathbf h}(f)$. The embedding for the norm
of a section gives
\[
 \|Q\|_{L^\infty(M,(T^*M)^{\otimes2}\otimes f^*TM;\mathbf h)}\leq C_q\sum_{j=0}^{q}\|\nabla^jQ\|_{L^2(M,(T^*M)^{\otimes(j+2)}\otimes f^*TM;\mathbf h)}.
\]
This estimate uses the induced metric connection: apply the scalar
Sobolev and Kato inequalities successively to
$Q,\nabla Q,\ldots,\nabla^{q-1}Q$. Its constants depend on the
geometry of the base. In particular, it does not require additional
bounds for the jets of $f$ itself.
The term $R$ in \eqref{eq:segunda-derivada-logaritmo-jets}
belongs to $L^2(M,(T^*M)^{\otimes2}\otimes TM;\mathbf h)$, since $X,\nabla X$ do and both are bounded.
Thus $X\in H^2(M,TM;\mathbf h)$.

Suppose that $X\in H^{j+1}(M,TM;\mathbf h)$, with $1\leq j\leq q$.
The transported connection differs from the fixed one by $a_X$,
given in \eqref{eq:variacion-conexion-transporte-curvatura}.
That formula controls $a_X$ in $H^j(M,T^*M\otimes\operatorname{End}TM;\mathbf h)\cap L^\infty(M,T^*M\otimes\operatorname{End}TM;\mathbf h)$ by
$X$ in $H^{j+1}(M,TM;\mathbf h)$ and its $C^1$ bound. On differentiating
\[
 \nabla(P_f^{-1}Q)=P_f^{-1}\nabla Q-a_X*(P_f^{-1}Q),
\]
the remaining terms are products of derivatives of $a_X$ and
$Q$ of total order at most $j-1$.
Hölder and the preceding estimates control them by the norms
of $X$ already obtained and the data of $Q$ up to order $q$.
This gives $P_f^{-1}Q\in H^j(M,(T^*M)^{\otimes2}\otimes TM;\mathbf h)$.
Also, $R(x,X,\nabla X)\in H^j(M,(T^*M)^{\otimes2}\otimes TM;\mathbf h)$ by the product rule.
Equation \eqref{eq:segunda-derivada-logaritmo-jets} gives
$X\in H^{j+2}(M,TM;\mathbf h)$. Each step produces a polynomial in the data
with no constant term: every summand retains a factor among
$X,\nabla X,Q$ or their derivatives. Induction up to $j=q$ proves
\eqref{eq:control-logaritmo-por-jets}.

The identities are first obtained for smooth maps.
Local $H^s(M,M)$ regularity allows approximation in each chart and
passage to weak derivatives. The uniform bounds are summed over
the cover to give the global norms.
Finally, starting from $X\in H^s(M,TM;\mathbf h)$ and differentiating the
exponential gives the converse by the same product bounds.
Subtracting the expression at $X=0$ ensures integrability
and convergence to zero.
\end{proof}

\begin{proposition}[Change of reference in the Sobolev component]
\label{prop:cambio-referencia-uniforme-rebanada}
Let $M$ (of dimension $m$) be smooth and without boundary, and let
$\mathbf g_0,\mathbf g$ be metrics of bounded geometry
of all orders. Fix $r\in\mathbb N$, with $r>m/2+2$, and assume
$\mathbf g\in\operatorname{Met}^{2,r}_{\mathbf g_0}(M)$.
For every fixed tensor bundle $\mathbf F\to M$ and $0\leq j\leq r$, the $H^j(M,\mathbf F)$
spaces of the two metrics agree, with equivalent norms.
Moreover,
\begin{equation}
\label{eq:identificacion-referencias-metricas-grupos}
 \operatorname{Met}^{2,r}_{\mathbf g_0}(M)
       =\operatorname{Met}^{2,r}_{\mathbf g}(M),
 \qquad
 \mathcal D^{r+1}_0(M,\mathbf g_0)
       =\mathcal D^{r+1}_0(M,\mathbf g).
\end{equation}
The first identification is a diffeomorphism of affine Hilbert
open sets; the second is an isomorphism of topological groups.
\end{proposition}

\begin{proof}
Write $\nabla_0=\nabla^{\mathbf g_0}$,
$\nabla=\nabla^{\mathbf g}$ and $A=\nabla-\nabla_0$.
Uniform positivity, the Sobolev embedding and identity
\eqref{eq:diferencia-levi-civita-metricas} give
\begin{equation}
\label{eq:cotas-relativas-cambio-referencia-metricas}
 C^{-1}\mathbf g_0\leq\mathbf g\leq C\mathbf g_0,
 \qquad A\in H^{r-1}(M,T^*M\otimes\operatorname{End}TM;\mathbf g_0)
\end{equation}
for some $C\geq1$. Indeed, the inverse metric is a multiplier
in the $H^r(M,S^2T^*M)$ regime, and $\nabla_0\mathbf g$ belongs to
$H^{r-1}(M,T^*M\otimes S^2T^*M;\mathbf g_0)$. Proposition~\ref{prop:comparacion-normas-metricas-componente}
proves norm equivalence up to order $r$.
In particular, $A\in H^{r-1}(M,T^*M\otimes\operatorname{End}TM;\mathbf g)$ as well.
If $\mathbf h=\mathbf g_0+b$, then
$\mathbf h-\mathbf g=b-(\mathbf g-\mathbf g_0)$ belongs to the
same $H^r(M,S^2T^*M)$ space with either norm. Pointwise comparison identifies
the uniform positivity conditions.
This proves the first equality; its affine transitions are
translations and a bicontinuous linear identity, hence smooth.

For the groups, set $s=r+1$ and compare the jets of a map,
without assuming equivalence of tensor norms at order $s$.
Let $f=\operatorname{Exp}^{\mathbf g_0}_{\operatorname{id}}X_0$,
with $X_0$ in a small ball of $H^s(M,TM;\mathbf g_0)$.
Choose its radius so that $f$ is a diffeomorphism,
$|df|_{\mathbf g_0}\leq2\sqrt m$, and both displacements lie in
the normal domains of the preceding lemma. This is possible because
$d_{\mathbf g}(x,f(x))\leq\sqrt C\,|X_0(x)|_{\mathbf g_0}$
and the embedding controls the $C^1$ norm.
Let $P_0,P$ be the radial transports for the two metrics.
Comparison of norms and volumes gives
\[
 \|d_{\mathbf g}(\,\cdot\,,f)\|_{L^2(M,\mathbf g)}
 \leq C^{(m+2)/4}\|X_0\|_{L^2(M,TM;\mathbf g_0)}.
\]

In a uniform normal chart for $\mathbf g_0$, the coefficients of
both connections are bounded: those of the second are those of
the first plus $A$, and $A\in C_b^1(M,T^*M\otimes\operatorname{End}TM;\mathbf g_0)$ by the embedding.
The radial $\mathbf g$ geodesic has length, measured with
$\mathbf g_0$, at most $C\,d_{\mathbf g_0}(x,f(x))$.
After reducing the radius, both geodesics remain in the same
uniform chart. The linear transport equation and Gronwall give
\[
 |P^{-1}P_0-I|_{\mathbf g_0}
       \leq C_2d_{\mathbf g_0}(x,f(x)).
\]
The constant depends on the connection bounds, the radius and
$C$, and is independent of $x$ and the map in the chosen ball.
Since $|df|$ is bounded,
\[
 \|K_{\mathbf g}(f)\|_{L^2(M,\operatorname{End}TM;\mathbf g)}
 \leq C_3\bigl(\|K_{\mathbf g_0}(f)\|_{L^2(M,\operatorname{End}TM;\mathbf g_0)}
                         +\|X_0\|_{L^2(M,TM;\mathbf g_0)}\bigr).
\]
These data tend to zero when $X_0\to0$ in $H^s(M,TM;\mathbf g_0)$.

The second derivative of the map satisfies the exact identity
\begin{equation}
\label{eq:cambio-referencia-segunda-derivada-mapeo}
 \begin{split}
 Q_{\mathbf g}(f)(U,V)={}&Q_{\mathbf g_0}(f)(U,V)\\
 &+A_{f(x)}(df_xU,df_xV)-df_x(A_x(U,V)),
 \qquad U,V\in T_xM.
 \end{split}
\end{equation}
Transport the target fibers by $P_0^{-1}$, and let $B_f$
denote the sum of the last two terms, transported in this way.
The composition and product estimates give
\[\begin{gathered}
 B_f\in H^{r-1}(M,(T^*M)^{\otimes2}\otimes TM;\mathbf g_0),
 \\ B_f\longrightarrow0\quad\text{in }H^{r-1}(M,(T^*M)^{\otimes2}\otimes TM;\mathbf g_0)
                              \\ \text{if }X_0\to0\text{ in }H^s(M,TM;\mathbf g_0).
\end{gathered}\]
In detail, $A\in H^{r-1}(M,T^*M\otimes\operatorname{End}TM;\mathbf g_0)$ and composition by $f$ is continuous
at that level, by the estimates in
Theorem~\ref{teo:eichhorn-difeomorfismos-abiertos}.
The factor $P_0^{-1}df=I+K_{\mathbf g_0}(f)$ is a multiplier
of $H^{r-1}(M,\mathbf F;\mathbf g_0)$, with integrable part in $H^r(M,\operatorname{End}TM;\mathbf g_0)$. Composition of
$A$, with its indices transported, and these factors are continuous
in the stated norms. At $f=\operatorname{id}$ the two terms
are $A$ and $-A$, so they cancel. Continuity for Sobolev $A$
follows by approximating it in $H^{r-1}(M,T^*M\otimes\operatorname{End}TM;\mathbf g_0)$ by smooth compactly
supported tensors and using the locally uniform bounds for the
composition operators.

We must also compare the connections used to differentiate
\eqref{eq:cambio-referencia-segunda-derivada-mapeo}.
After transport by $P_0$, the difference from the fixed
$\mathbf g_0$ connection contains the transport difference
$a_{X_0}$, the actions of $A$ on the two base indices, and
the action of $A_{f(x)}df_x$ on the target index.
Proposition~\ref{prop:transporte-homotopia-conexion-sobolev}
gives $a_{X_0}\in H^r(M,T^*M\otimes\operatorname{End}TM;\mathbf g_0)$; the other coefficients belong to
$H^{r-1}(M,T^*M\otimes\operatorname{End}TM;\mathbf g_0)$ by the preceding paragraph. Their norms are bounded
on the chosen ball. Expansion
\eqref{eq:expansion-cambio-conexion-tensor-metricas}, applied
to these transported bundles, and the module lemma give
\[\begin{aligned}
 \sum_{j=0}^{r-1}\|\nabla^jQ_{\mathbf g}(f)\|_{L^2(M,(T^*M)^{\otimes(j+2)}\otimes f^*TM;\mathbf g)}\\
 &\leq C_4\|P_0^{-1}Q_{\mathbf g_0}(f)\|_{H^{r-1}(M,(T^*M)^{\otimes2}\otimes TM;\mathbf g_0)}
\\
 &\quad+C_4\|B_f\|_{H^{r-1}(M,(T^*M)^{\otimes2}\otimes TM;\mathbf g_0)}.
\end{aligned}\]
Indeed, derivatives of the coefficients only reach order $r-2$;
the products are controlled from the fixed connection, using
$r-1>m/2$. Bounded geometry is not assumed for the connection
induced by the Sobolev map.
The first term on the right tends to zero by the converse
part of the lemma, applied to $\mathbf g_0$, and the second by
continuity of $B_f$.

We have proved that $J_{s,\mathbf g}(f)\to0$ when
$X_0\to0$ in $H^s(M,TM;\mathbf g_0)$.
Lemma~\ref{lem:carta-exponencial-control-jets-referencia} implies
\[
 X:=\log^{\mathbf g}_{\operatorname{id}}f\in H^s(M,TM;\mathbf g),
 \qquad \|X\|_{H^s(M,TM;\mathbf g)}\longrightarrow0.
\]
Local membership in $H^s(M,M)$ is preserved by changes of smooth
charts on any compact set, so the local hypothesis of the
lemma holds. Thus a neighborhood of the identity in the first
group lies in the second, and the identification is continuous
at the identity. Interchanging the metrics gives the reverse
inclusion and continuity.

The intersection of the two groups contains a neighborhood
of the identity in each. It is therefore an open and closed
subgroup of each identity component; connectedness forces it
to be the whole component. This proves the second equality
in \eqref{eq:identificacion-referencias-metricas-grupos}.
Continuity at the identity of a homomorphism of topological
groups implies continuity everywhere by translation.
The identification is therefore a homeomorphism.
\end{proof}

The change of reference for the group corresponds to
Proposition~2.22 of \cite[pp.~94--95]{EichhornSlice}.
The comparison uses $A\in H^{r-1}(M,T^*M\otimes\operatorname{End}TM;\mathbf g_0)$ and the jets of the map up
to order $r+1$. The group identification is asserted in the
topological category. If uniform bounds for all derivatives
of $A$ are also imposed, expansion
\eqref{eq:expansion-cambio-conexion-tensor-metricas} gives
equivalence of $H^j(M,\mathbf F)$ norms on each tensor bundle $\mathbf F\to M$ for every $j\in\mathbb N_0$.

\begin{corollary}[Slice after a change of reference]
\label{cor:rebanada-metrica-cambio-referencia-uniforme}
Under the hypotheses of
Proposition~\ref{prop:cambio-referencia-uniforme-rebanada}, suppose
that $m\geq2$ and $K_{\mathbf g}\leq-c$ for some $c>0$.
There are a smooth Hilbert submanifold $\mathcal S$ containing
$\mathbf g$ and an invariant open set $\mathcal U$ in
$\operatorname{Met}^{2,r}_{\mathbf g_0}(M)$ such that
\[
 \mathcal S\times\mathcal D^{r+1}_0(M,\mathbf g_0)
       \longrightarrow\mathcal U,
 \qquad (\mathbf s,f)\longmapsto f^*\mathbf s,
\]
is a homeomorphism. The action on $\mathcal U$ is proper, and each
orbit contained in $\mathcal U$ is closed in
$\operatorname{Met}^{2,r}_{\mathbf g_0}(M)$.
The orbit of $\mathbf g$ is a smooth Hilbert submanifold
with complemented tangent space.
\end{corollary}

\begin{proof}
Apply Theorem~\ref{teo:rebanada-armonica-metrica-referencia}
with reference $\mathbf g$.
The first identification in
\eqref{eq:identificacion-referencias-metricas-grupos} sends its
invariant open set and its slice to the same sets of
metrics with reference $\mathbf g_0$.
The second identification preserves the group elements,
their topology, and pullback. Consequently, the equivariant
product and its inverse remain continuous.
Properness of the action and closedness of each orbit are
preserved by these homeomorphisms. The smooth affine
identification also preserves the submanifold structures
of the slice and the orbit of $\mathbf g$.
\end{proof}

To assemble these tubes, fix a smooth Riemannian manifold
$(M,\mathbf g_0)$ (of dimension $m$) without boundary and of bounded geometry of all orders,
with $m\geq2$ and $K_{\mathbf g_0}\leq-c_0$ for some $c_0>0$.
Let $r\in\mathbb N$, $r>m/2+2$, and write
\[
 \mathscr M^r:=\operatorname{Met}^{2,r}_{\mathbf g_0}(M),
 \qquad \mathscr G:=\mathcal D^{r+1}_0(M,\mathbf g_0).
\]
The Sobolev norms of tensors use $\mathbf g_0$ as reference.
Define
\begin{equation}
\label{eq:region-negativa-sobolev-y-centros-suaves}
\begin{split}
 \Omega_r^-:=\{\mathbf h\in\mathscr M^r\mid {}&
 \text{there is }c_{\mathbf h}>0\text{ such that }
 K_{\mathbf h}(x,\sigma)\leq-c_{\mathbf h}\\
 &\text{for every }x\in M,\quad
                  \sigma\in\operatorname{Gr}_2(T_xM)\},\\
 \mathscr C_\infty^-:=\{\mathbf h\in\Omega_r^-\mid {}&
 \mathbf h\text{ is smooth and of bounded geometry of all orders}\}.
\end{split}
\end{equation}
The embedding $H^r(M,S^2T^*M)\hookrightarrow C_b^2(M,S^2T^*M)$ allows the curvature
of every $\mathbf h\in\mathscr M^r$ to be interpreted classically.

\begin{lemma}[Openness and approximation of negatively curved metrics]
\label{lem:aproximacion-metricas-negativas-centros-suaves}
The set $\Omega_r^-$ is open and invariant under $\mathscr G$.
The set $\mathscr C_\infty^-$ is dense in $\Omega_r^-$.
More precisely, every $\mathbf h\in\Omega_r^-$ is an $H^r(M,S^2T^*M)$ limit
of metrics in $\mathscr C_\infty^-$ that agree with $\mathbf g_0$
outside a compact set.
\end{lemma}

\begin{proof}
Fix $\mathbf h\in\Omega_r^-$ and $0<\alpha\leq1$ such that
$\alpha\mathbf g_0\leq\mathbf h\leq\alpha^{-1}\mathbf g_0$.
Let $C_S\geq1$ be an embedding constant satisfying
\[
 \|\mathbf a\|_{C_b^2(M,S^2T^*M;\mathbf g_0)}
 :=\max_{0\leq j\leq2}\sup_{x\in M}
             |(\nabla^{\mathbf g_0})^j\mathbf a(x)|_{\mathbf g_0}
 \leq C_S\|\mathbf a\|_{H^r(M,S^2T^*M;\mathbf g_0)}.
\]
The curvature formula in terms of the metric, its inverse and
its first two derivatives gives a constant $C_{\mathbf h}\geq1$
such that
\[
 \sup_{\substack{x\in M\\\sigma\in\operatorname{Gr}_2(T_xM)}}
 |K_{\mathbf k}(x,\sigma)-K_{\mathbf h}(x,\sigma)|
 \leq C_{\mathbf h}\|\mathbf k-\mathbf h\|_{C_b^2(M,S^2T^*M;\mathbf g_0)}
\]
if $\|\mathbf k-\mathbf h\|_{C_b^2(M,S^2T^*M;\mathbf g_0)}
<\displaystyle\min\{1,\alpha/2\}$. Indeed, in the uniform charts the
coefficients and their first two derivatives remain bounded,
the metric $\mathbf k$ is at least $(\alpha/2)\mathbf g_0$, and
the inverse metrics are uniformly bounded. To compute sectional
curvature, choose a $\mathbf g_0$-orthonormal basis of the plane;
the denominator in its expression for $\mathbf k$ is at least
$\alpha^2/4$. These bounds give the stated constant by the
product rule and the identity for the difference of inverses.

Choose
\[
 \varepsilon_{\mathbf h}:=C_S^{-1}
       \min\{1,\alpha/2,c_{\mathbf h}/(2C_{\mathbf h})\}.
\]
If $\|\mathbf k-\mathbf h\|_{H^r(M,S^2T^*M;\mathbf g_0)}
<\varepsilon_{\mathbf h}$, then $\mathbf k$ is uniformly
positive and $K_{\mathbf k}\leq-c_{\mathbf h}/2$.
This proves openness. Invariance follows from the naturality of
curvature under pullback and the continuity of the action in
Proposition~\ref{prop:distancia-sobolev-invariante-metricas}.

Write $\mathbf h=\mathbf g_0+\mathbf a$.
The cutoffs in \eqref{eq:cola-densidad-Wkp-geometria-acotada-cap13}
and Theorem~\ref{meyers-serrin-haz} provide tensors
$\mathbf a_j\in C_c^\infty(M,S^2T^*M)$ converging to
$\mathbf a$ in $H^r(M,S^2T^*M)$. For all sufficiently large $j$,
$\mathbf h_j:=\mathbf g_0+\mathbf a_j$ belongs to the ball above.
Thus $\mathbf h_j\in\Omega_r^-$ and is complete by its uniform
equivalence to $\mathbf g_0$.

Each $\mathbf h_j$ has bounded geometry of all orders.
To check the injectivity radius, set
$K_j=\operatorname{supp}\mathbf a_j$ and choose
$0<\rho<\operatorname{inj}(M,\mathbf g_0)/4$.
If $K_j$ is nonempty, its closed neighborhood
\[
 L_j:=\{x\in M\mid d_{\mathbf g_0}(x,K_j)\leq2\rho\}
\]
is compact by completeness. For $x\notin L_j$, the
$\mathbf g_0$-geodesics starting at $x$ with length less than
$\rho$ remain outside $K_j$, where they are also
$\mathbf h_j$-geodesics. Their exponential maps agree on that
ball, so $\operatorname{inj}_{\mathbf h_j}(x)\geq\rho$.
On $L_j$, compactness gives a positive lower bound for the
injectivity radius of the smooth metric $\mathbf h_j$.
The curvature derivatives are bounded on $L_j$ by smoothness
and agree with those of $\mathbf g_0$ outside $K_j$.
If $K_j$ is empty, $\mathbf h_j=\mathbf g_0$.
Hence $\mathbf h_j\in\mathscr C_\infty^-$, and convergence proves
density.
\end{proof}

\begin{lemma}[Estimates near the identity]
\label{lem:estimaciones-identidad-propiedad-propia}
With the fixed reference $\mathbf g_0$ and integer $r$, write
$\nabla^0=\nabla^{\mathbf g_0}$. There are $\delta>0$ and $C>0$
with the following property. Let $X\in H^{r+1}(M,TM)$ satisfy
$\|X\|_{C_b^1(M,TM;\mathbf g_0)}<\delta$, and put
$f_X(x)=\exp_x^{\mathbf g_0}X(x)$. Suppose that $f_X$ is a
diffeomorphism, and denote radial parallel transport by $P_X$.
For tensors $A,B\in H^{r-1}(M,T^*M\otimes T^*M\otimes TM)$,
define
\[
 \begin{aligned}
 A[X](U,V)&:=P_X^{-1}A_{f_X(x)}(df_XU,df_XV),\\
 B\langle X\rangle(U,V)&:=P_X^{-1}df_X\bigl(B_x(U,V)\bigr),
 \qquad U,V\in T_xM.
 \end{aligned}
\]
Then
\begin{align}
 \|A[X]\|_{H^{r-1}(M,T^*M\otimes\operatorname{End}TM;\mathbf g_0)}
 &\leq C\|A\|_{H^{r-1}(M,T^*M\otimes\operatorname{End}TM;\mathbf g_0)}
\nonumber\\
 &\quad+C\|A\|_{C_b^1(M,T^*M\otimes\operatorname{End}TM;\mathbf g_0)}\|X\|_{H^r(M,TM;\mathbf g_0)},
 \label{eq:composicion-lineal-orden-superior-propia}\\
 \|B\langle X\rangle\|_{H^{r-1}(M,T^*M\otimes\operatorname{End}TM;\mathbf g_0)}
 &\leq C\|B\|_{H^{r-1}(M,T^*M\otimes\operatorname{End}TM;\mathbf g_0)}
\nonumber\\
 &\quad+C\|B\|_{C_b^0(M,T^*M\otimes\operatorname{End}TM;\mathbf g_0)}\|X\|_{H^r(M,TM;\mathbf g_0)}.
 \label{eq:producto-lineal-orden-superior-propia}
\end{align}
Moreover, if
$S_X:=P_X^{-1}\nabla^0df_X$ and $\|S_X\|_{C_b^0(M,(T^*M)^{\otimes2}\otimes TM;\mathbf g_0)}\leq R$, there is
$C_R>0$, independent of $X$, such that
\begin{equation}
\label{eq:recuperacion-orden-superior-difeomorfismo}
 \|X\|_{H^{r+1}(M,TM;\mathbf g_0)}
 \leq C_R\bigl(\|S_X\|_{H^{r-1}(M,(T^*M)^{\otimes2}\otimes TM;\mathbf g_0)}+\|X\|_{H^r(M,TM;\mathbf g_0)}\bigr).
\end{equation}
All these norms are taken over $M$ with reference $\mathbf g_0$.
\end{lemma}

\begin{proof}
Reduce $\delta$ so that each pair $x,f_X(x)$ and its radial geodesic
lie in a larger uniform chart. The jet bounds of
Lemma~\ref{lem:jets-mixtos-exponencial-orden-finito} give common bounds
for $df_X$, its inverse, and their Jacobians. In these charts,
$f_X(x)=x+v(x)$, and the Sobolev norms of $v$ are controlled by those
of $X$, through superposition with uniform coefficients vanishing
at $X=0$.

We detail the composition estimate. Let $s=r-1>m/2+1$.
Differentiating $a\circ f_X$ a total of $s$ times gives contractions
of terms of the form
\[
 (D^\ell a)\circ f_X
       \prod_{\nu=1}^{\ell}D^{j_\nu}f_X,
 \qquad
 1\leq\ell\leq s,\quad j_\nu\geq1,
 \quad\sum_{\nu=1}^{\ell}j_\nu=s.
\]
For $\ell>1$, interpolation between $Da\in L^\infty(\mathbb R^m)$ and
$D^sa\in L^2(\mathbb R^m)$ assigns the first factor the exponent
$2(s-1)/(\ell-1)$. For each $j_\nu>1$, interpolate between
$Df_X\in L^\infty(\mathbb R^m)$ and $D^sv\in L^2(\mathbb R^m)$ with exponent
$2(s-1)/(j_\nu-1)$. First-order factors are taken in
$L^\infty(\mathbb R^m)$. The reciprocals of these exponents sum to $1/2$, since
\[
 (\ell-1)+\sum_{\nu=1}^{\ell}(j_\nu-1)=s-1.
\]
Hölder, change of variables, and Young's inequality give a bound
by a constant times
$\|a\|_{H^s(\mathbb R^m)}+\|a\|_{C^1}\|v\|_{H^s(\mathbb R^m)}$.
The endpoint cases $\ell=1,s$ are estimated directly with
$Da$ or $Df_X$ in $L^\infty(\mathbb R^m)$. Chart cutoffs introduce only lower
derivatives. Summing over the uniform cover, whose multiplicity is
bounded, gives the same global estimate. The additional $df_X$
factors are controlled by the product inequality in the proof of
Lemma~\ref{lem:carta-exponencial-control-jets-referencia}; the
highest-order term uses $\|X\|_{H^{s+1}(M,TM;\mathbf g_0)}$. Radial transport is a
superposition coefficient with the same jet bounds. This proves
\eqref{eq:composicion-lineal-orden-superior-propia}.
The second estimate follows directly from the product with
$P_X^{-1}df_X=I+K_X$, where
$\|K_X\|_{H^{r-1}(M,\operatorname{End}TM;\mathbf g_0)}\leq C\|X\|_{H^r(M,TM;\mathbf g_0)}$ and
$\|K_X\|_{C_b^0(M,\operatorname{End}TM;\mathbf g_0)}\leq C\delta$.

For the last bound, use
\eqref{eq:segunda-derivada-logaritmo-jets}, with reference
$\mathbf g_0$:
\[
 (\nabla^0)^2X=A_X^{-1}S_X+R_0(x,X,\nabla^0X).
\]
The coefficient $A_X^{-1}$ has uniform jets. The expression
$R_0$ has degree at most two in $\nabla^0X$, with smooth uniform
coefficients, and $R_0(x,0,0)=0$. The same product and
superposition estimates give
\[
 \begin{aligned}
 \|A_X^{-1}S_X\|_{H^{r-1}(M,(T^*M)^{\otimes2}\otimes TM;\mathbf g_0)}
   &\leq C\bigl(\|S_X\|_{H^{r-1}(M,(T^*M)^{\otimes2}\otimes TM;\mathbf g_0)}+R\|X\|_{H^{r-1}(M,TM;\mathbf g_0)}\bigr),\\
 \|R_0(x,X,\nabla^0X)\|_{H^{r-1}(M,(T^*M)^{\otimes2}\otimes TM;\mathbf g_0)}
   &\leq C\|X\|_{H^r(M,TM;\mathbf g_0)}.
 \end{aligned}
\]
Together with orders zero and one of $X$, these inequalities prove
\eqref{eq:recuperacion-orden-superior-difeomorfismo}.
\end{proof}

\begin{proposition}[Convergence of the pullback and the diffeomorphism]
\label{prop:rigidez-sucesiones-casi-isometricas}
Let $\mathbf h\in\Omega_r^-$ and let $(f_n)_{n\in\mathbb N}$ be a
sequence in $\mathscr G$. If $f_n^*\mathbf h\to\mathbf h$ in
$H^r(M,S^2T^*M;\mathbf g_0)$, then $f_n\to\operatorname{id}_M$ in $\mathscr G$.
Here $\mathbf h$ has the $H^r(M,S^2T^*M;\mathbf g_0)$ regularity of $\mathscr M^r$ and,
by Sobolev embedding, is of class $C^2$.
\end{proposition}

\begin{proof}
Put $\mathbf k_n=f_n^*\mathbf h$ and choose $c>0$ with
$K_{\mathbf h}\leq-c$. Uniform equivalence to $\mathbf g_0$
makes $\mathbf h$ complete. Convergence in $C_b^2(M,S^2T^*M;\mathbf g_0)$ allows us to assume
\[
 \tfrac12\mathbf h\leq\mathbf k_n\leq2\mathbf h.
\]
Since $f_n:(M,\mathbf k_n)\to(M,\mathbf h)$ is an isometry,
\begin{equation}
\label{eq:tension-casi-isometria-metrica-sobolev}
 \tau_{\mathbf h,\mathbf h}(f_n)
   =df_n\bigl(\operatorname{tr}_{\mathbf h}
                (\nabla^{\mathbf k_n}-\nabla^{\mathbf h})\bigr),
 \qquad
 t_n:=\|\tau_{\mathbf h,\mathbf h}(f_n)\|_{C_b^0(M,f_n^*TM;\mathbf h)}\longrightarrow0.
\end{equation}
Indeed, $\|df_n\|_{\mathbf h\to\mathbf h}\leq\sqrt2$, and the
connection-difference formula bounds the right-hand side by
$C_{\mathbf h}\|\mathbf k_n-\mathbf h\|_{C_b^1(M,S^2T^*M;\mathbf g_0)}$.

Lift $f_n$ from the identity to the universal covering space of
each component of $M$. The argument of
Lemma~\ref{lem:unicidad-armonica-componente-sobolev} shows that the
lifted displacement descends to a bounded function on $M$ tending
to zero outside compact sets. In that argument paths are constructed
using $\mathbf g_0$; uniform equivalence allows their lengths to be
measured using $\mathbf h$. Thus the function
\[
 u_n(x)=\tfrac12 d_{\widetilde{\mathbf h}}
                         (z,\widetilde f_n(z))^2,
 \qquad \pi(z)=x,
\]
attains a maximum if it takes a positive value. Cartan--Hadamard
and second variation apply to this $C^2$ metric: the connection
is $C^1$, the Jacobi equations have continuous coefficients, and
squared distance is $C^2$; the second-variation formula contains
only derivatives of the metric up to order two.

Let $x$ be a point of positive maximum, and let
$d=d_{\widetilde{\mathbf h}}(z,\widetilde f_n(z))>0$.
Parametrize the geodesic between these points on $[0,1]$.
Its velocity $T$ has norm $d$. For an
$\mathbf h$-orthonormal basis $e_1,\ldots,e_m$ at $x$, let $J_i$
be the Jacobi fields obtained by varying both endpoints through
the identity and $f_n$. Then $J_i(0)=e_i$.
The sum of the boundary terms is
$\langle T(1),\tau_{\mathbf h,\mathbf h}(f_n)\rangle$.
If $Y_i$ is the perpendicular part of $J_i$, the curvature bound gives
\[
 \int_0^1\bigl(|D_tJ_i|^2
       -\langle R^{\mathbf h}(J_i,T)T,J_i\rangle\bigr)dt
 \geq\int_0^1\bigl(|D_tY_i|^2+cd^2|Y_i|^2\bigr)dt.
\]
Write $\lambda=\sqrt c\,d$ and
$b(t)=\lambda\tanh(\lambda(1-t))$. In a parallel frame, the
identity $b^2-b'=\lambda^2$, with $b(1)=0$, gives
\[
 \int_0^1\bigl(|Y'|^2+\lambda^2|Y|^2\bigr)dt
 =\int_0^1|Y'+bY|^2dt+\lambda\tanh(\lambda)|Y(0)|^2.
\]
Since $\displaystyle\sum_{i=1}^{m}|Y_i(0)|^2=m-1$, second
variation at the maximum implies
\begin{equation}
\label{eq:control-desplazamiento-por-tension}
 0\geq(m-1)\sqrt c\,d\tanh(\sqrt c\,d)-d\,t_n.
\end{equation}
For every $\varepsilon>0$, if
$t_n<(m-1)\sqrt c\tanh(\sqrt c\,\varepsilon)$, this inequality
prevents the maximum displacement from being at least $\varepsilon$.
Thus the lifted displacement, and also
$\displaystyle\sup_{x\in M}d_{\mathbf g_0}(x,f_n(x))$, tend to zero.

In larger uniform charts, the connection identity for an isometry
reads, for $1\leq a,i,j\leq m$,
\[
 \partial_i\partial_j f_n^a
 =\sum_{b=1}^{m}\Gamma(\mathbf k_n)^b_{ij}\partial_b f_n^a
  -\sum_{b,c=1}^{m}\Gamma(\mathbf h)^a_{bc}(f_n)
                              \partial_i f_n^b\partial_j f_n^c.
\]
The bounds on $df_n$, the metrics, and their first derivatives
give a uniform bound on $D^2f_n$ in these charts. If
$v_n=f_n-\operatorname{id}$ in coordinates, Taylor's formula on
segments of length $\ell$ inside the larger charts gives, on the
smaller charts,
\[
 \|Dv_n\|_{C^0}\leq C\bigl(\ell^{-1}\|v_n\|_{C^0}+\ell\bigr),
 \qquad 0<\ell<\rho_0/2,
\]
where $\rho_0>0$ is a common radius. First choose $\ell$ small
and then $n$ large. This gives uniform convergence of $df_n$ to
the differential of the identity. For large $n$, we may write
\[
 f_n=f_{X_n},\qquad X_n(x)=\log_x^{\mathbf g_0}f_n(x),
 \qquad \|X_n\|_{C_b^1(M,TM;\mathbf g_0)}\longrightarrow0.
\]
The regularity $X_n\in H^{r+1}(M,TM)$ follows from the mapping-space
charts and Sobolev superposition. It can also be checked using a
finite product of exponential charts representing $f_n$: outside
a compact set all its fields are small in $C^1$. For their first
derivatives, uniform vanishing outside compact sets follows from
the same disjoint-ball argument applied to $\nabla^0X_j\in L^2(M,T^*M\otimes TM;\mathbf g_0)$,
whose derivative is bounded. The logarithm of the composition is
then a uniform superposition vanishing when those fields are zero.
On the remaining compact set, local Sobolev regularity suffices.

The Korn identity of
Lemma~\ref{lem:coercividad-operadores-rebanada-metrica} holds for
$\mathbf h$ of class $C^2$ by the same integration-by-parts
calculation and density. Since
$\operatorname{Ric}_{\mathbf h}\leq-(m-1)c\mathbf h$, it gives a
constant $a_{\mathbf h}>0$ such that
\begin{equation}
\label{eq:korn-c2-metrica-sobolev}
 a_{\mathbf h}\|X\|_{H^1(M,TM;\mathbf g_0)}
 \leq\|\mathcal L_X\mathbf h\|_{L^2(M,S^2T^*M;\mathbf g_0)}.
\end{equation}
Here we used equivalence of the first-order norms and of the
zeroth-order tensor norms for the two metrics.
Taylor's formula with integral remainder in the coordinate
expression for pullback gives $C_{\mathbf h}>0$ and
$\delta_{\mathbf h}>0$ such that
\begin{equation}
\label{eq:resto-pullback-c1-h1}
 \|f_X^*\mathbf h-\mathbf h-\mathcal L_X\mathbf h\|_{L^2(M,S^2T^*M;\mathbf g_0)}
 \leq C_{\mathbf h}\|X\|_{C_b^1(M,TM;\mathbf g_0)}\|X\|_{H^1(M,TM;\mathbf g_0)},
 \qquad \|X\|_{C_b^1(M,TM;\mathbf g_0)}<\delta_{\mathbf h}.
\end{equation}
Indeed, each remainder term is bounded pointwise by
$C_{\mathbf h}(|X|+|\nabla^0X|)^2$: only the first two derivatives
of $\mathbf h$ and the uniform exponential-map jets enter.
Taking the magnitude of one factor in $L^\infty(M)$ and that of the other in $L^2(M)$ proves the bound.
For sufficiently large $n$,
$C_{\mathbf h}\|X_n\|_{C_b^1(M,TM;\mathbf g_0)}<a_{\mathbf h}/2$, and therefore
\[
 \|X_n\|_{H^1(M,TM;\mathbf g_0)}
 \leq\frac{2}{a_{\mathbf h}}\|\mathbf k_n-\mathbf h\|_{L^2(M,S^2T^*M;\mathbf g_0)}
 \longrightarrow0.
\]

We now recover the highest order. Put
$A=\nabla^{\mathbf h}-\nabla^0$ and
$B_n=\nabla^{\mathbf k_n}-\nabla^0$.
The connection-difference formula gives
$B_n\to A$ in $H^{r-1}(M,T^*M\otimes\operatorname{End}TM;\mathbf g_0)$ and, by embedding, in $C_b^1(M,T^*M\otimes\operatorname{End}TM;\mathbf g_0)$.
The preceding connection identity takes the intrinsic form
\begin{equation}
\label{eq:conexion-isometria-control-sobolev}
 S_{X_n}=B_n\langle X_n\rangle-A[X_n].
\end{equation}
Its $C_b^0(M,(T^*M)^{\otimes2}\otimes TM;\mathbf g_0)$ norm is bounded by a common constant $R$.
Lemma~\ref{lem:estimaciones-identidad-propiedad-propia} gives a
constant $C_0$, independent of $n$, such that
\[
 \|X_n\|_{H^{r+1}(M,TM;\mathbf g_0)}\leq C_0(1+\|X_n\|_{H^r(M,TM;\mathbf g_0)}).
\]
Interpolation on the uniform cover provides, for every
$\eta>0$, a constant $C_\eta$ such that
\[
 \|X_n\|_{H^r(M,TM;\mathbf g_0)}\leq\eta\|X_n\|_{H^{r+1}(M,TM;\mathbf g_0)}
                              +C_\eta\|X_n\|_{H^1(M,TM;\mathbf g_0)}.
\]
Choosing $C_0\eta\leq1/2$ gives a uniform bound in
$H^{r+1}(M,TM;\mathbf g_0)$. The multiplicative interpolation inequality then gives
\[
 \|X_n\|_{H^r(M,TM;\mathbf g_0)}
 \leq C\|X_n\|_{H^1(M,TM;\mathbf g_0)}^{1/r}
              \|X_n\|_{H^{r+1}(M,TM;\mathbf g_0)}^{(r-1)/r}\longrightarrow0.
\]
Continuity of composition at the lower level, applied with
$r>m/2+1$ and tensors of order $H^{r-1}(M,T^*M\otimes\operatorname{End}TM;\mathbf g_0)$, implies
$A[X_n]\to A$ and $B_n\langle X_n\rangle\to A$ in $H^{r-1}(M,T^*M\otimes\operatorname{End}TM;\mathbf g_0)$.
On the uniform cover, this continuity follows as in
Lemma~\ref{lem:composicion-sobolev-autocontenida}: approximate the
tensor by smooth compactly supported tensors and control the errors
using the composition norms, which are uniformly bounded because
$X_n\to0$ in $H^r(M,TM;\mathbf g_0)$.
Thus \eqref{eq:conexion-isometria-control-sobolev} gives
$S_{X_n}\to0$ in $H^{r-1}(M,(T^*M)^{\otimes2}\otimes TM;\mathbf g_0)$. Estimate
\eqref{eq:recuperacion-orden-superior-difeomorfismo} finally proves
$X_n\to0$ in $H^{r+1}(M,TM;\mathbf g_0)$, as asserted.
\end{proof}

\begin{theorem}[Proper action on the entire negatively curved region]
\label{teo:accion-propia-region-negativa-completa}
The action of $\mathscr G$ on $\Omega_r^-$ is free and proper.
For every $\mathbf h\in\Omega_r^-$, its orbit is closed in
$\mathscr M^r$, and the orbit map
$f\mapsto f^*\mathbf h$ is a homeomorphism from $\mathscr G$
onto that orbit.
\end{theorem}

\begin{proof}
The preceding proposition, applied to a constant sequence, shows
that $f^*\mathbf h=\mathbf h$ implies
$f=\operatorname{id}_M$. It also proves continuity of the inverse
orbit map at $\mathbf h$: otherwise, there would be a neighborhood
$V$ of the identity and a sequence $f_n\notin V$ with
$d_r(f_n^*\mathbf h,\mathbf h)<1/n$.
Translations in the group and the metric space give continuity
at every other point of the orbit.

We prove closedness. Suppose that $f_n^*\mathbf h\to\mathbf k$
in $\mathscr M^r$. Invariance of $d_r$ gives
\[
 d_r((f_m\circ f_n^{-1})^*\mathbf h,\mathbf h)
       =d_r(f_m^*\mathbf h,f_n^*\mathbf h)\longrightarrow0
                         \quad(n,m\to\infty).
\]
By the inverse continuity just proved,
$f_m\circ f_n^{-1}\to\operatorname{id}_M$, uniformly for
sufficiently large $m,n$, in the group topology.

We detail the resulting convergence. Choose $\rho>0$ so that the
ball of radius $2\rho$ in $H^{r+1}(M,TM)$ lies in an exponential
chart of diffeomorphisms. On a smaller ball, composition and the
integral formula for superposition differences give a constant
$C_\rho$ such that
\begin{equation}
\label{eq:cauchy-derecha-coordenadas-difeomorfismos}
 \begin{gathered}
 \|Y'-Y\|_{H^{r+1}(M,TM;\mathbf g_0)}\leq C_\rho\|Z\|_{H^{r+1}(M,TM;\mathbf g_0)}\quad\text{if}\quad
 f_{Y'}=f_Z\circ f_Y,\\
 \|Y\|_{H^{r+1}(M,TM;\mathbf g_0)},\|Y'\|_{H^{r+1}(M,TM;\mathbf g_0)},\|Z\|_{H^{r+1}(M,TM;\mathbf g_0)}<\rho.
 \end{gathered}
\end{equation}
To justify uniformity, composition by $f_Y$ is a bounded operator
on $H^{r+1}(M,TM;\mathbf g_0)$ with a common norm bound on that ball: the Jacobians,
the bi-Lipschitz bounds, and the $H^{r+1}(M,TM;\mathbf g_0)$ norm of $Y$ are controlled.
In exponential coordinates, the difference between
$f_Z\circ f_Y$ and $f_Y$ is a smooth superposition vanishing at
$Z=0$. Its differential with respect to $Z\circ f_Y$, together
with the bound on $Y$, has a uniform norm bound. Integrating along
the segment from zero to $Z\circ f_Y$ proves
\eqref{eq:cauchy-derecha-coordenadas-difeomorfismos}.

Fix a sufficiently large $N$ and write
$f_n\circ f_N^{-1}=f_{Y_n}$, with
$\|Y_n\|_{H^{r+1}(M,TM;\mathbf g_0)}<\rho/2$, for $n\geq N$.
For large $m,n$, also write
$f_m\circ f_n^{-1}=f_{Z_{mn}}$, where
$\|Z_{mn}\|_{H^{r+1}(M,TM;\mathbf g_0)}\to0$ as $m,n\to\infty$.
Since $f_{Y_m}=f_{Z_{mn}}\circ f_{Y_n}$, the preceding estimate
makes $(Y_n)$ Cauchy in $H^{r+1}(M,TM;\mathbf g_0)$.
Its limit $Y$ has norm at most $\rho/2$, so it remains in the chart.
Thus $f_n\to f_Y\circ f_N=:f\in\mathscr G$.
Continuity of the action gives $\mathbf k=f^*\mathbf h$.
This proves closedness in $\mathscr M^r$ and convergence of the
diffeomorphism for every convergent sequence in that orbit.

Now consider
\[
 \Psi:\Omega_r^-\times\mathscr G\longrightarrow
                    \Omega_r^-\times\Omega_r^- ,
 \qquad \Psi(\mathbf h,f)=(\mathbf h,f^*\mathbf h).
\]
Let $C$ be compact in the target and let
$(\mathbf h_n,f_n)_{n\in\mathbb N}$ be a sequence in
$\Psi^{-1}(C)$. Passing to a subsequence, we may assume
\[
 (\mathbf h_n,f_n^*\mathbf h_n)\longrightarrow(\mathbf h,\mathbf k)
                                                        \quad\text{in }C.
\]
Invariance of $d_r$ implies
\begin{equation}
\label{eq:convergencia-orbita-fija-accion-propia}
 d_r(f_n^*\mathbf h,\mathbf k)
 \leq d_r(\mathbf h,\mathbf h_n)
                +d_r(f_n^*\mathbf h_n,\mathbf k)\longrightarrow0.
\end{equation}
What we have proved gives $f_n\to f\in\mathscr G$ and
$\mathbf k=f^*\mathbf h$. Consequently,
$(\mathbf h_n,f_n)\to(\mathbf h,f)\in\Psi^{-1}(C)$.
The domain is metrizable, so sequential compactness of
$\Psi^{-1}(C)$ is equivalent to compactness. Hence $\Psi$ is proper.
\end{proof}

For each $\mathbf g\in\mathscr C_\infty^-$, choose a tube from
Corollary~\ref{cor:rebanada-metrica-cambio-referencia-uniforme},
with invariant open set $\mathcal U_{\mathbf g}\subseteq\Omega_r^-$
and slice $\mathcal S_{\mathbf g}$ containing $\mathbf g$.
The choice within $\Omega_r^-$ is possible by
Theorem~\ref{teo:rebanada-armonica-metrica-referencia}.
Set
\begin{equation}
\label{eq:union-tubos-metricas-negativas}
 \mathscr U_r^-:=\bigcup_{\mathbf g\in\mathscr C_\infty^-}
                         \mathcal U_{\mathbf g}.
\end{equation}

\begin{theorem}[Proper action on the union of tubes]
\label{teo:accion-propia-union-tubos-metricas}
The set $\mathscr U_r^-$ is an invariant open dense subset of
$\Omega_r^-$ containing $\mathscr C_\infty^-$.
The action of $\mathscr G$ on $\mathscr U_r^-$ is free and proper.
Each of its orbits is closed in $\mathscr M^r$, and the orbit map
$f\mapsto f^*\mathbf h$ is a homeomorphism from $\mathscr G$
onto the orbit of $\mathbf h$.
\end{theorem}

\begin{proof}
Openness and invariance are preserved under unions.
Each center belongs to its tube, so
$\mathscr C_\infty^-\subseteq\mathscr U_r^-$; density follows
from Lemma~\ref{lem:aproximacion-metricas-negativas-centros-suaves}.
Theorem~\ref{teo:accion-propia-region-negativa-completa} gives
triviality of isotropy, closedness of orbits, and the assertion
about their parametrizations. Restricting the action to this
invariant open set preserves properness: the preimage of a compact
subset of $\mathscr U_r^-\times\mathscr U_r^-$ equals its preimage
for the action on $\Omega_r^-$ and is contained in
$\mathscr U_r^-\times\mathscr G$.
\end{proof}

\begin{corollary}[Global distance and quotient charts]
\label{cor:cociente-topologico-union-tubos-metricas}
Let $q:\Omega_r^-\to\Omega_r^-/\mathscr G$ be the projection.
The formula
\begin{equation}
\label{eq:distancia-cociente-metricas-negativas}
 \overline d_r([\mathbf h],[\mathbf k])
 :=\inf_{f\in\mathscr G}d_r(\mathbf h,f^*\mathbf k)
\end{equation}
defines a distance inducing the quotient topology.
The space $\Omega_r^-/\mathscr G$ is Hausdorff and second countable.
The dense open subset $q(\mathscr U_r^-)$ is a topological Hilbert
manifold, and $q|_{\mathscr U_r^-}$ is a topological principal bundle
with structure group $\mathscr G$.
For every center $\mathbf g\in\mathscr C_\infty^-$, the restriction
of $q$ to $\mathcal S_{\mathbf g}$ is a homeomorphism onto the
open set $q(\mathcal U_{\mathbf g})$.
\end{corollary}

\begin{proof}
Invariance of $d_r$ shows that the infimum is independent of the
representatives. Indeed, replacing $\mathbf h$ by $a^*\mathbf h$
and $\mathbf k$ by $b^*\mathbf k$ gives
\[
 d_r(a^*\mathbf h,f^*b^*\mathbf k)
   =d_r(\mathbf h,(b\circ f\circ a^{-1})^*\mathbf k),
 \qquad a,b,f\in\mathscr G,
\]
and $b\circ f\circ a^{-1}$ ranges over all of $\mathscr G$.
Symmetry follows by replacing $f$ with $f^{-1}$.
For $\mathbf h,\mathbf k,\boldsymbol\ell\in\Omega_r^-$ and
$a,b\in\mathscr G$, the triangle inequality gives
\[
 \overline d_r([\mathbf h],[\boldsymbol\ell])
 \leq d_r(\mathbf h,(b\circ a)^*\boldsymbol\ell)
 \leq d_r(\mathbf h,a^*\mathbf k)
                         +d_r(\mathbf k,b^*\boldsymbol\ell).
\]
Taking the infima over $a,b\in\mathscr G$ gives the triangle
inequality for $\overline d_r$.
If $\overline d_r([\mathbf h],[\mathbf k])=0$, then for every
$n\in\mathbb N$ there is $f_n\in\mathscr G$ such that
$d_r(\mathbf h,f_n^*\mathbf k)<1/n$. Closedness of the orbit
of $\mathbf k$ implies $[\mathbf h]=[\mathbf k]$.
The converse follows by choosing an element that takes
$\mathbf k$ to $\mathbf h$.

The projection $q$ is open: for an open set
$V\subseteq\Omega_r^-$,
\[
 q^{-1}(q(V))=\bigcup_{f\in\mathscr G}f^*V
\]
is open. Moreover, for $\mathbf h\in\Omega_r^-$ and
$\varepsilon>0$,
\begin{equation}
\label{eq:bolas-distancia-cociente-metricas}
 q\bigl(B_{d_r}(\mathbf h,\varepsilon)\cap\Omega_r^-\bigr)
   =B_{\overline d_r}([\mathbf h],\varepsilon).
\end{equation}
One inclusion uses the identity element; the other uses the fact
that an infimum less than $\varepsilon$ has some term less than
$\varepsilon$. These equalities and continuity of $q$ for
$\overline d_r$, which follows from
$\overline d_r([\mathbf h],[\mathbf k])\leq d_r(\mathbf h,\mathbf k)$,
prove equality of the topologies.

In the tube coordinates
$\mathcal S_{\mathbf g}\times\mathscr G\to\mathcal U_{\mathbf g}$,
the orbits are exactly the sets
$\{\mathbf s\}\times\mathscr G$. Their quotient is
$\mathcal S_{\mathbf g}$ with its own topology, since projection
of the product onto the first factor is open.
Invariance of $\mathcal U_{\mathbf g}$ and openness of $q$
identify this quotient with $q(\mathcal U_{\mathbf g})$.
Charts of the Hilbert submanifold $\mathcal S_{\mathbf g}$ thus
give topological charts of $q(\mathscr U_r^-)$. Equivariance of the
product, where the action is $(\mathbf s,a)\cdot f
=(\mathbf s,a\circ f)$, provides the principal trivializations
of $q|_{\mathscr U_r^-}$.

Openness of $q$ and density of $\mathscr U_r^-$ in $\Omega_r^-$
show that $q(\mathscr U_r^-)$ is open and dense.
The distance $\overline d_r$ ensures the Hausdorff property.
The space $H^r(M,S^2T^*M;\mathbf g_0)$ is separable; hence
$\Omega_r^-$ has a countable base of open sets.
Their images under the open map $q$ form a countable base of
the quotient. Finally, all points of $\Omega_r^-$ have the
same orbit type, determined by the stabilizer
$\{\operatorname{id}_M\}$.
\end{proof}

\section{Diffeomorphism groups, Hilbert scales, and ILH}
\label{sec:grupos-difeomorfismos-ilh}

\subsection{The Sobolev group at fixed regularity}

Let $M$ be a smooth compact manifold of dimension $m$ with or without boundary and let
$s>\frac{m}{2}+1$. If $\partial M=\varnothing$, the set
\[
\operatorname{Diff}^s(M):=
\{\eta\in H^s(M,M)\mid \eta\text{ is bijective and }
\eta^{-1}\in H^s(M,M)\}
\]
is open in $H^s(M,M)$ and is therefore a Hilbert manifold. If
$\partial M\neq\varnothing$, introduce the double $DM$, which is a
finite-dimensional manifold without boundary. In this case,
$\operatorname{Diff}^s(M)$ consists of the $C^1$ diffeomorphisms of $M$
whose coordinate representations and those of their inverses are $H^s(M,M)$; in
boundary charts, use the restriction spaces on the half-space.
Work near one such diffeomorphism
$\eta\in\operatorname{Diff}^s(M)$ within $H^s(M,DM)$.
Choose on $DM$ a product metric in a collar of $\partial M$; thus
$\partial M$ is totally geodesic. In the exponential chart centered at
$\eta$, the condition of mapping boundary to boundary is equivalent to
\[
 \operatorname{Tr}(\mathbf X)\in
 H^{s-\frac12}(\partial M,(\eta|_{\partial M})^*T\partial M).
\]
Indeed, a sufficiently short geodesic with both endpoints on
$\partial M$ is the geodesic for the induced metric on that hypersurface;
its initial vector is tangent. The reverse implication follows because the
hypersurface is totally geodesic.

The operator taking the normal component of the trace,
\[
 H^s(M,\eta^*TDM)\longrightarrow
 H^{s-\frac12}(\partial M,(\eta|_{\partial M})^*N\partial M),
\]
is continuous and surjective. The trace extension operator in the collar,
followed by inclusion of the normal bundle, provides a continuous
right inverse. Its kernel is therefore a closed complemented subspace.
Sobolev trivializations of the pullback bundle are justified using
$H^{s-1/2}(\partial M)$, which is an algebra in this range.

It remains to check that maps in this chart take values in $M$ and are
diffeomorphisms. In collar coordinates with inward coordinate $t\geq0$,
the normal component of $\eta$ vanishes at $t=0$, and its normal derivative is
strictly positive there. Compactness of the boundary gives a positive
lower bound on a smaller collar. A sufficiently small perturbation
in $C^1$ preserves this bound and, by integration from $t=0$, maps the
interior collar into the interior. Outside the collar, the image under $\eta$ lies
at positive distance from the boundary, a condition preserved by
$C^0$ closeness.

Specify a $C^1$ neighborhood that guarantees injectivity.
Fix a compatible distance $d$ on $M$ and choose finitely many charts
$\phi_i\colon U_i\to B_i$, where $B_i$ is a convex ball or
convex half-ball, and open sets $V_i$ with
$\overline{V_i}\subset U_i$ covering $M$.
Also choose target charts $\psi_i$ containing
$\eta(\overline{U_i})$. If $a_i\in B_i$ is the center and
$\eta_i=\psi_i\circ\eta\circ\phi_i^{-1}$, set
$A_i=D\eta_i(a_i)$. Shrinking the charts gives
\[
 \sup_{u\in B_i}\|A_i^{-1}D\eta_i(u)-I\|\leq\tfrac14.
\]
The compact sets $\eta(\overline{U_i})$ have positive margins
inside the target charts. Since there are finitely many charts, one
$C^1$ neighborhood of $\eta$ ensures, for every perturbation
$\zeta$, that $\zeta_i=\psi_i\circ\zeta\circ\phi_i^{-1}$ is
defined on $B_i$ and satisfies
$\sup_{B_i}\|A_i^{-1}D\zeta_i-I\|\leq1/2$ for every $i$.
Integration along the segment joining $u,v\in B_i$ gives
\[
 A_i^{-1}\bigl(\zeta_i(u)-\zeta_i(v)\bigr)
 =(u-v)+\int_0^1
 \bigl(A_i^{-1}D\zeta_i(v+t(u-v))-I\bigr)(u-v)\,dt,
\]
and therefore
\begin{equation}\label{eq:inyectividad-uniforme-perturbacion-compacta}
 \bigl|A_i^{-1}(\zeta_i(u)-\zeta_i(v))\bigr|
 \geq\tfrac12|u-v|.
\end{equation}
Thus $\zeta$ is injective on each $U_i$, and its differential is
invertible there.

Let $\lambda>0$ be a Lebesgue number of the cover $\{V_i\}$.
If $d(p,q)<\lambda$, both points lie in one $U_i$.
For the remaining pairs, the set
\[
 P_\lambda=\{(p,q)\in M\times M\mid d(p,q)\geq\lambda\}
\]
is compact. If it is nonempty, injectivity of $\eta$ gives
$a:=\min_{P_\lambda}d(\eta(p),\eta(q))>0$.
Shrink the neighborhood to ensure
$\sup_p d(\zeta(p),\eta(p))<a/3$; then
$d(\zeta(p),\zeta(q))\geq a/3$ on $P_\lambda$.
This proves global injectivity.

The image is open by the inverse function theorem; at the boundary,
use preservation of the boundary and the positive normal derivative
already proved. It is closed by compactness.
Since $M$ has finitely many components, choose one point in each
and shrink the $C^0$ neighborhood so that its image under $\zeta$
stays in the component of its image under $\eta$.
The open and closed image of $\zeta$ then meets every component,
so it equals $M$. Regularity of the inverse in $H^s(M,M)$ and its
continuity are proved in
Lemma~\ref{lem:inversion-sobolev-real-autocontenida}.

The embedding $H^s(M,M)\hookrightarrow C^1(M,M)$ makes this argument applicable in a
Sobolev neighborhood of $\eta$. We obtain a chart of
$\operatorname{Diff}^s(M)$ on the kernel of the normal trace. Its tangent
space at $\eta$ consists of sections of $\eta^*TM$ whose trace is tangent to
the boundary. This construction is Lemma~6.6 of
\cite{EbinMarsden1970}; it uses a Hilbert manifold without boundary modeled
on that kernel.

Localization near diffeomorphisms is part of the argument. The set
of all maps of pairs $M\to M$ may have additional constraints
imposed by positivity of the normal coordinate. For example, at the
constant zero map from $[0,1]$ to itself, every differentiable variation with
a parameter of both signs has zero velocity at every point, since the value
zero is a minimum. The preceding linear model is not obtained there.

Inversion is studied together with composition in
\cite{EbinMarsden1970,InciKappelerTopalov2012}. We prove the required
regularity and continuity here, including real orders and the local
treatment of the boundary.

\begin{lemma}[Regularity and continuity of Sobolev inversion]
\label{lem:inversion-sobolev-real-autocontenida}
Let $M$ be a smooth compact manifold of dimension $m$, with or without
boundary, and let $s>m/2+1$ be real. If $\eta\colon M\to M$ is a $C^1$
diffeomorphism whose local representations belong to $H^s(M,M)$, then the
representations of $\eta^{-1}$ also belong to $H^s(M,M)$. Moreover, if
$\eta_n\to\eta$ in $H^s(M,M)$ and all these maps are $C^1$ diffeomorphisms,
then $\eta_n^{-1}\to\eta^{-1}$ in $H^s(M,M)$.
\end{lemma}

\begin{proof}
The point that cannot be assumed is the regularity of the inverse.
The change of variables estimate of order zero and the fractional
estimate of order less than one only require a bi-Lipschitz map,
not an inverse already known to be $H^s(M,M)$. We begin with these two
estimates.

Let $U,V$ be bounded coordinate domains and let $\psi\colon V\to U$
be a $C^1$ diffeomorphism onto its image which is bi-Lipschitz. If its
Lipschitz constants and Jacobian bounds are fixed, then, for
$0\leq\sigma<1$,
\begin{equation}\label{eq:pullback-bajo-orden-inversion}
 \|w\circ\psi\|_{H^\sigma(V)}\leq C_\sigma\|w\|_{H^\sigma(U)}.
\end{equation}
For $\sigma=0$ this is change of variables in the integral of
$|w\circ\psi|^2$. For $0<\sigma<1$, use the equivalent norm defined
by the difference seminorm. If $x=\psi(y)$ and $x'=\psi(y')$, then
$|y-y'|\geq |x-x'|/\operatorname{Lip}(\psi)$, and hence
\[
 \int_V\int_V
 \frac{|w(\psi(y))-w(\psi(y'))|^2}{|y-y'|^{m+2\sigma}}\,dy\,dy'
 \leq C\int_U\int_U
 \frac{|w(x)-w(x')|^2}{|x-x'|^{m+2\sigma}}\,dx\,dx'.
\]
Restriction norms on smooth domains or half-balls are compared with
these norms by the local extension operators already constructed.

Put $\psi=\eta^{-1}$, initially only of class $C^1$.
Choose smaller charts with closure in larger charts, so that $\psi$
maps the closure of the smaller target chart into the relative interior
of a source chart. At the boundary use half-ball charts. After shrinking
them, the representations of $\eta$ and $\psi$ are bi-Lipschitz with
finite constants. Finitely many such pairs cover $M$.

In one of the source charts define
\[
 A_\eta(x)=(D\eta(x))^{-1}.
\]
Since $s-1>m/2$, smooth superposition on
$H^{s-1}(U)$, established in
Corollary~\ref{cor:superposicion-hilbert-orden-real}, gives
$A_\eta\in H^{s-1}_{\mathrm{loc}}(U,\mathbb R^{m\times m})$ up to the boundary.
Indeed, $D\eta$ has compact image in the open set of invertible matrices
over each smaller chart. Matrix inversion can be replaced there by a
smooth function on the whole matrix space, using a cutoff equal to one
on a neighborhood of that image. This is not inversion on a Sobolev
mapping space, but only inversion of finite matrices.

Write $s=k+\sigma$, with $k\in\mathbb N$ and $0\leq\sigma<1$.
For a word $\alpha=(i_1,\ldots,i_q)$ of indices in
$\{1,\ldots,m\}$ we construct a coordinate field $B_\alpha$ such that
\begin{equation}\label{eq:recurrencia-derivadas-inversa-sobolev}
 \partial_{i_q}\cdots\partial_{i_1}\psi=B_\alpha\circ\psi,
 \qquad B_\alpha\in H^{s-q}_{\mathrm{loc}}(U,\mathbb R^m),\qquad 1\leq q\leq k.
\end{equation}
For $q=1$, the identity $D\psi=A_\eta\circ\psi$ gives
$B_{(i)}=A_\eta e_i$. If $B_\alpha$ has been constructed and $q<k$, set
\[
 B_{(\alpha,i)}
 =\displaystyle\sum_{j=1}^{m}(\partial_jB_\alpha)(A_\eta)_{ji}.
\]
This formula is used in smaller charts; to estimate the norms we insert
cutoffs equal to one there. The cutoffs are fixed, and their products
and derivatives are bounded operators on the spaces in use.
Lemma~\ref{lem:modulo-hilbert-orden-real-composicion} gives
\[
 H^{s-q-1}(U)\cdot H^{s-1}(U)\subseteq H^{s-q-1}(U),
 \qquad 0\leq s-q-1\leq s-1,
\]
which proves the required regularity of $B_{(\alpha,i)}$.

We must also justify the weak derivative identity. When $q<k$,
$B_\alpha$ belongs locally to $H^1(U,\mathbb R^m)$. Approximation by smooth functions
in $H^1(U,\mathbb R^m)$, the $L^2(U)$ change of variables bound and boundedness of $D\psi$
allow passage to the limit in the classical chain rule. This gives
$\partial_i(B_\alpha\circ\psi)=B_{(\alpha,i)}\circ\psi$
in distributions. Starting from the first-order $C^1$ identity proves
\eqref{eq:recurrencia-derivadas-inversa-sobolev} up to order $k$, without
presupposing higher derivatives of $\psi$.

For $q\leq k$, change of variables shows that these derivatives are
in $L^2(V,\mathbb R^m)$. For $q=k$ we also have $B_\alpha\in H^\sigma(U,\mathbb R^m)$, and
\eqref{eq:pullback-bajo-orden-inversion} gives
$B_\alpha\circ\psi\in H^\sigma(V,\mathbb R^m)$. The characterization of $H^{k+\sigma}(V,\mathbb R^m)$
by derivatives of order $k$ and their fractional seminorms proves
$\psi\in H^s(V,\mathbb R^m)$ on each smaller chart. If $\sigma=0$, the $L^2(V,\mathbb R^m)$ conclusion
at order $k$ suffices. Finite localization proves the assertion on $M$.

We turn to continuity. The embedding $H^s(M,M)\hookrightarrow C^1(M,M)$ gives
$\eta_n\to\eta$ in $C^1$. First $\psi_n:=\eta_n^{-1}\to\psi$
uniformly. Indeed, if $x_n=\psi_n(y)$, then
\[
 d(\eta(x_n),y)\leq\displaystyle\sup_{x\in M}d(\eta(x),\eta_n(x));
\]
uniform continuity of $\eta^{-1}$ turns this into a uniform bound for
$d(x_n,\psi(y))$ tending to zero. The identity
$D\psi_n=(D\eta_n\circ\psi_n)^{-1}$, applied in the chosen charts,
then gives $C^1$ convergence. In particular, the same smaller charts,
Jacobian bounds and bi-Lipschitz constants can be used for every
sufficiently large $n$.

Continuity of matrix superposition and of the module product, applied
to the preceding recurrence, gives
\[
 B_{\alpha,n}\longrightarrow B_\alpha
       \quad\text{in }H^{s-q}(U,\mathbb R^m),\qquad |\alpha|=q\leq k,
\]
with the same local cutoffs. For every $w\in H^\sigma(U)$ we have
$w\circ\psi_n\to w\circ\psi$ in $H^\sigma(V)$, $0\leq\sigma<1$:
approximate $w$ by a smooth function $w_\varepsilon$ in $H^\sigma(U,\mathbb R^m)$.
The uniform bound \eqref{eq:pullback-bajo-orden-inversion} controls the
two errors involving $w-w_\varepsilon$. For fixed $w_\varepsilon$,
$C^1$ convergence of $\psi_n$ gives $C^1$ convergence of the compositions,
and hence $H^\sigma(V,\mathbb R^m)$ convergence on the smaller bounded domain.
First let $n\to\infty$ and then $\varepsilon\to0$. Thus
\[
 B_{\alpha,n}\circ\psi_n\longrightarrow B_\alpha\circ\psi
 \quad\text{in }H^\sigma(V,\mathbb R^m)\quad\text{if }|\alpha|=k.
\]
The same argument with $\sigma=0$ treats derivatives of lower order.
Together with uniform convergence of $\psi_n$, these convergences prove
$\psi_n\to\psi$ in $H^s(M,M)$. The proof in boundary charts is identical:
approximations are obtained by extending the coordinate functions, but
changes of variables are performed only between the interior domains.
We do not extend the global diffeomorphism to another one whose $H^s(M,M)$
regularity has been assumed in advance.
\end{proof}

\begin{theorem}[Ebin--Marsden]
\label{teo:ebin-marsden-diff-sobolev}
Let $M$ be a smooth compact manifold with or without boundary. If
$s>\frac{m}{2}+1$, the preceding space
is a topological group. For
$\eta\in\operatorname{Diff}^s(M)$, right translation
$R_\eta(\varphi)=\varphi\circ\eta$ is smooth. Left translation
and inversion, however, are not generally smooth at the same level. More
precisely, for $j\in\mathbb N_0$,
\[
L_\eta\colon \operatorname{Diff}^s(M)\longrightarrow\operatorname{Diff}^s(M)
\quad\text{is }C^j\text{ if }\eta\in\operatorname{Diff}^{s+j}(M),
\]
and inversion is $C^j$ as a map
$\operatorname{Diff}^{s+j}(M)\longrightarrow\operatorname{Diff}^s(M)$.
Thus $\operatorname{Diff}^s(M)$ is not a Hilbert Lie group in the
usual sense.
\end{theorem}

\begin{proof}
For $j=0$, continuity of composition is
Lemma~\ref{lem:composicion-sobolev-autocontenida}, while continuity of
inversion is
Lemma~\ref{lem:inversion-sobolev-real-autocontenida}. Thus
$\operatorname{Diff}^s(M)$ is a topological group. Right translation by a
fixed diffeomorphism is smooth: in mapping-space charts it is pullback by a
fixed diffeomorphism and its derivatives are continuous linear operators of
the same type; changes of charts only add the smooth superposition operators
already established.

Let $j\geq1$. Lemma~\ref{lem:composicion-sobolev-perdida-autocontenida}
gives the $C^j$ regularity of
\[
 H^{s+j}(M,M)\times\operatorname{Diff}^s(M)
 \longrightarrow H^s(M,M),\qquad (\eta,\varphi)\longmapsto\eta\circ\varphi,
\]
in the relevant local charts. Hence, for
$\eta\in\operatorname{Diff}^{s+j}(M)$, left translation $L_\eta$ is $C^j$.

For inversion we do not differentiate an identity before knowing that the
inverse is a map between the stated spaces. This has already been proved in
Lemma~\ref{lem:inversion-sobolev-real-autocontenida}. Near
$(\eta,\eta^{-1})$ consider
\[
 \mathcal F(\varphi,\psi)=\varphi\circ\psi-\operatorname{id}.
\]
As a map from a neighborhood in
$\operatorname{Diff}^{s+j}(M)\times\operatorname{Diff}^{s}(M)$ to $H^s(M,TM)$, it
is $C^j$. Its derivative in the second variable is
\[
 D_2\mathcal F(\eta,\eta^{-1})\mathbf X
 = (d\eta\circ\eta^{-1})\mathbf X.
\]
The coefficient $d\eta\circ\eta^{-1}$ belongs to $H^s(M,\operatorname{End}TM)$ and takes values in
invertible matrices; its matrix inverse also belongs to $H^s(M,\operatorname{End}TM)$, by smooth
matrix superposition and the inversion lemma just proved. Since $H^s(M)$ is an
algebra, multiplication by this coefficient is a bounded isomorphism of the
$H^s(M,TM)$ tangent space. The same statement holds in boundary charts on the
subspace preserving tangential trace. The Banach implicit function theorem
therefore produces a unique local $C^j$ map
$\varphi\mapsto\psi(\varphi)$ such that
$\varphi\circ\psi(\varphi)=\operatorname{id}$. Set-theoretic uniqueness of
the inverse gives $\psi(\varphi)=\varphi^{-1}$. Consequently,
\[
 \operatorname{inv}\colon\operatorname{Diff}^{s+j}(M)
 \longrightarrow\operatorname{Diff}^{s}(M)
\]
is $C^j$. Only after this regularity has been established may one
differentiate the identity and obtain
\begin{equation}\label{eq:diferencial-inversion-sobolev-justificado}
 D(\operatorname{inv})(\eta)\mathbf X
 =-\bigl(d\eta\circ\eta^{-1}\bigr)^{-1}
   \bigl(\mathbf X\circ\eta^{-1}\bigr).
\end{equation}
The derivative loss prevents these operations, in general, from defining a
Hilbert Lie group at one fixed level.
\end{proof}

\begin{example}[Derivative loss already occurs on the circle]
\label{ej:inversion-composicion-no-c1-nivel-fijo}
Let $M=S^1$ and $s>3/2$. Choose a real function
$u\in H^s(S^1)\setminus H^{s+1}(S^1)$ with sufficiently small $C^1$ norm;
for instance, take a Fourier series having exactly $H^s(S^1)$ regularity and
multiply it by a small constant. Then
\[
 \eta(x)=x+u(x)\pmod{2\pi}
\]
is an $H^s(S^1,S^1)$ diffeomorphism. If left translation
$L_\eta\colon\operatorname{Diff}^s(S^1)\to\operatorname{Diff}^s(S^1)$
were differentiable at the identity, using the translation curve
$\tau_t(x)=x+t$ would force
\[
 \frac{\eta\circ\tau_t-\eta}{t}
\]
to converge in $H^s(S^1)$. In $H^{s-1}(S^1)$ it converges, by the fundamental theorem
of calculus, to $\eta'$. Uniqueness of the limit under the continuous
inclusion $H^s(S^1)\hookrightarrow H^{s-1}(S^1)$ would imply $\eta'\in H^s(S^1)$, hence
$\eta\in H^{s+1}(S^1,S^1)$, a contradiction.

Inversion is not generally $C^1$ at this fixed level either. Indeed,
Lemma~\ref{lem:inversion-sobolev-real-autocontenida}, also applied at order
$s+1$, shows that if $\eta^{-1}$ belonged to $H^{s+1}(S^1,S^1)$ then $\eta$ would
belong to $H^{s+1}(S^1,S^1)$. Thus $\eta^{-1}\notin H^{s+1}(S^1,S^1)$. For
$\eta_t=\tau_t\circ\eta$ one has
$\eta_t^{-1}=\eta^{-1}\circ\tau_{-t}$, and the same difference-quotient
argument would force $(\eta^{-1})'\in H^s(S^1)$ if inversion were $C^1$ at
$\eta$. This example explains why the regularity recovered by moving up the
scale is structural rather than a defect of the proof.
\end{example}

\begin{lemma}[Composition at lower levels and differentiation of curves]
\label{lem:composicion-niveles-inferiores-curvas}
Let $M$ be a smooth compact manifold of dimension $m$, with or without boundary,
and let $s>m/2+1$. For $0\leq r\leq s$, composition defines a continuous
map
\[
 H^r(M,\mathbb R^d)\times\operatorname{Diff}^s(M)
 \longrightarrow H^r(M,\mathbb R^d),\qquad(F,\eta)\longmapsto F\circ\eta.
\]
For each $\eta_0$, there exist a neighborhood $\mathcal U$ and a constant
$C_{r,\mathcal U}$ with
\begin{equation}
\label{eq:composicion-nivel-inferior-cota}
 \|F\circ\eta\|_{H^r(M,\mathbb R^d)}
 \leq C_{r,\mathcal U}\|F\|_{H^r(M,\mathbb R^d)},\qquad\eta\in\mathcal U.
\end{equation}
If $F\in C^1(I,H^s(M,\mathbb R^d))$ and
$\eta\in C^1(I,\operatorname{Diff}^s(M))$, then
$z=F\circ\eta^{-1}$ belongs to
$C(I,H^s(M,\mathbb R^d))\cap C^1(I,H^{s-1}(M,\mathbb R^d))$ and
\begin{equation}
\label{eq:derivada-curva-composicion-inversa}
 \partial_tz=(\partial_tF)\circ\eta^{-1}-dz(\mathbf v),
 \qquad\mathbf v=\dot\eta\circ\eta^{-1}.
\end{equation}
In local frames, the same conclusions hold for sections along
$\eta$; using covariant derivatives, the formula is
$\partial_t\mathbf z=(\nabla_t\mathbf F)\circ\eta^{-1}
-\nabla_{\mathbf v}\mathbf z$.
\end{lemma}

\begin{proof}
The bound \eqref{eq:composicion-nivel-inferior-cota} and local boundedness of its constants are exactly Lemma~\ref{lem:composicion-sobolev-autocontenida}. No continuity of inversion is needed here to bound Jacobians.

To check joint continuity, let $F_j\to F$ in
$H^r(M,\mathbb R^d)$ and $\eta_j\to\eta$ in
$\operatorname{Diff}^s(M)$. Choose $G\in C^\infty(M,\mathbb R^d)$
with $\|F-G\|_{H^r(M,\mathbb R^d)}<\varepsilon$. For large $j$,
all composition operators under consideration have a common bound,
and
\begin{align*}
 \|F_j\circ\eta_j-F\circ\eta\|_{H^r(M,\mathbb R^d)}
 &\leq C\|F_j-F\|_{H^r(M,\mathbb R^d)}+2C\varepsilon\\
 &\quad+\|G\circ\eta_j-G\circ\eta\|_{H^r(M,\mathbb R^d)}.
\end{align*}
The last term tends to zero even in $H^s(M,\mathbb R^d)$,
by smooth superposition with $G$. First let $j\to\infty$ and
then $\varepsilon\to0$.

For curves, continuity already proved gives
$z\in C(I,H^s(M,\mathbb R^d))$. Evaluation and embedding into
$C^1$ allow pointwise differentiation of the identity
$F(t,x)=z(t,\eta(t,x))$, giving
\eqref{eq:derivada-curva-composicion-inversa} pointwise.
Its right-hand side is a continuous curve in
$H^{s-1}(M,\mathbb R^d)$: the first term is continuous by composition,
and the second by the product and differentiation estimate in
Lemma~\ref{lem:calculo-sobolev-compacto-aplicaciones}.
Denote this side by $b(t)$. The integral
$z(t_*)+\displaystyle \int_{t_*}^tb(\tau)\,d\tau$ exists in
$H^{s-1}(M,\mathbb R^d)$. Point evaluation is continuous in this
space because $s-1>m/2$; evaluating gives the same integral identity
satisfied by $z(t,x)$. The two functions agree, and by the
fundamental theorem of calculus for the Bochner integral, $z$ is
of class $C^1$ in $H^{s-1}(M,\mathbb R^d)$ with derivative $b$.
For sections, perform the calculation in a finite cover of frames.
The connection coefficients add the smooth contractions of the
covariant derivative and satisfy the same product estimates;
the formulas agree on overlaps by the transformation rule
for a connection.
\end{proof}

\begin{theorem}[Sobolev diffeomorphisms in bounded geometry without boundary]
\label{teo:eichhorn-difeomorfismos-abiertos}
Let $(M,\mathbf{g})$ be a smooth oriented noncompact Riemannian manifold of dimension $m$
without boundary and of bounded geometry of order $k$. For
$1<p<\infty$ and $k>r>\frac{m}{p}+1$, let
\[
 \mathcal D^{p,r}(M)
 :=\left\{\boldsymbol{\eta}\in\Omega^{p,r}(M,M)\ \middle|\
 \begin{array}{l}
 \boldsymbol{\eta}\text{ is a }C^1\text{ orientation-preserving diffeomorphism},\\
 \displaystyle\inf_{x\in M}\lambda_{\min}(d\boldsymbol{\eta}_x)>0
 \end{array}\right\},
\]
where $\lambda_{\min}(d\boldsymbol{\eta}_x)
:=\displaystyle\inf_{\substack{v\in T_xM\\|v|_{\mathbf{g},x}=1}}|d\boldsymbol{\eta}_xv|_{\mathbf{g},\boldsymbol{\eta}(x)}$. Then
$\mathcal D^{p,r}(M)$ is open in $\Omega^{p,r}(M,M)$, and each component is
a Banach manifold of class $C^{k+1-r}$, a Hilbert manifold for $p=2$. The
identity component $\mathcal D^{p,r}_0(M)$ is a metrizable topological
group under composition and inversion.
\end{theorem}

Openness can be seen in the constructed charts. If $\eta$ belongs to the
preceding set and $h$ is sufficiently close in a chart,
$C^1$ closeness preserves a positive lower bound for $dh$. The uniform
bound on $d\eta^{-1}$ also implies that
$F=\eta^{-1}\circ h$ has small uniform displacement from the
identity. A map with bounded displacement on a complete
finite-dimensional manifold is proper: the preimage of a compact set lies in a
closed bounded metric neighborhood of that compact set, which is compact by
Hopf--Rinow. Since $F$ is a local diffeomorphism, it is a proper covering.
If the displacement is less than the injectivity radius, the homotopy
$F_t(x)=\exp_x(t\log_xF(x))$ joins it to the identity. On each component,
$F_*$ therefore induces an isomorphism of fundamental groups, and the covering
has one sheet. This gives injectivity and
surjectivity of $F$ and hence of $h$. The covering
and continuation arguments are those in the chapter on infinite-dimensional calculus.

We detail the estimates needed for the operations. Fix a neighborhood
of $\eta$ in the Sobolev charts. The embedding into $C_b^1(M,TM)$ and the
lower bound for $d\eta$ give common bounds for
$\|d\eta\|_\infty$, $\|(d\eta)^{-1}\|_\infty$, and the Jacobians
of $\eta$ and $\eta^{-1}$. For localization, use balls of a common
radius mapping into charts of a common radius. The multiplicity of
the image cover is bounded by the bi-Lipschitz bounds and the ball
volume bounds. This allows the local estimates to be summed. Norms on
$\mathbb R^m$ in this argument apply to expressions localized
by smooth cutoff functions compactly supported in the charts
and extended by zero outside their coordinate domains.

In a chart write $\eta=\operatorname{Exp}_f\mathbf X$ with smooth
center. For $2\leq j\leq r$, its coordinate derivatives of order $j$
split into a uniformly bounded part coming from the center and an
$L^{p_j}(\mathbb R^m)$ part, where
\[
 \frac1{p_j}=\frac{j-1}{(r-1)p}.
\]
Indeed, subtract the expression at $\mathbf X=0$ and use Taylor's
formula. Each remaining term contains a derivative of $\mathbf X$;
the embeddings and H\"older's inequality, with $r-1>m/p$, give the
stated bound. The norms are controlled by a polynomial in
$\|\mathbf X\|_{W^{r,p}(M,TM)}$ and the derivative bounds for $f$ through
order $r$.

For $u\in W^{r,p}(M,\mathbb R^d)$, first apply Faà di Bruno's
formula from Theorem~\ref{teo:faa-di-bruno-particiones} to smooth
coordinate representatives. Let $U_0,V_0\subseteq\mathbb R^m$
be the corresponding coordinate domains, with $\eta(U_0)\subseteq V_0$.
For $x\in U_0$, $v=(v_1,\ldots,v_j)\in(\mathbb R^m)^j$, and
$1\leq j\leq r$, we have
\[
\begin{split}
 D^j(u\circ\eta)(x)[v_1,\ldots,v_j]
 =\sum_{\substack{\pi\in\Pi_j\\\pi=\{B_1,\ldots,B_a\}}}
 D^au(\eta(x))
 [D_{B_1}\eta(x)[v],\ldots,D_{B_a}\eta(x)[v]].
\end{split}
\]
Here $\Pi_j$ is the set of partitions of $\{1,\ldots,j\}$ into
nonempty blocks, enumerated so that $\min B_1<\cdots<\min B_a$.
For any smooth coordinate representative $\zeta$ and
$B=\{i_1<\cdots<i_b\}$, we use
\[
 D_B\zeta(z)[v]:=D^b\zeta(z)[v_{i_1},\ldots,v_{i_b}].
\]
Symmetry of the differentials makes the evaluation independent
of the chosen order of the blocks. For a partition with $a$ blocks,
the orders satisfy
\[
 (a-1)+\sum_{B\in\pi}(|B|-1)=j-1\leq r-1.
\]
This is the count behind H\"older's inequality: apply the embedding
of $W^{r-1,p}(\mathbb R^m)$ to $Du$ and to the integrable parts of $D\eta$,
assigning each factor an exponent whose fraction of $1/p$ is at least
its derivative order divided by $r-1$. Choose these fractions to sum
to one. Bounded factors are taken in $L^\infty(\mathbb R^m)$.
Change of variables controls composition of each factor coming from
$u$; for the order zero term, explicitly,
\[
 \int_M|u(\eta(x))|^p\,d\lambda_{\mathbf g}(x)
 =\int_M|u(y)|^p J_{\eta^{-1}}(y)\,d\lambda_{\mathbf g}(y).
\]
Here $J_{\eta^{-1}}$ is the absolute Jacobian for the indicated
volume. Summing the localized estimates gives $C>0$, independent
of $u$ and uniform on the fixed neighborhood, such that
\[
 \|u\circ\eta\|_{W^{r,p}(M,\mathbb R^d)}
 \leq C\|u\|_{W^{r,p}(M,\mathbb R^d)}.
\]
For sections the identities are read in local
frames, adding the connection coefficients already controlled.

This calculation does not assume Sobolev regularity of inversion.
If $\chi=\eta^{-1}$ is the $C^1$ inverse, the identity
$\eta\circ\chi=\operatorname{id}$ first gives
$D\chi(y)=(D\eta(\chi(y)))^{-1}$. If $\eta$ is smooth,
for $y\in\eta(U_0)$, $v_1,\ldots,v_j\in\mathbb R^m$, and
$2\leq j\leq r$, applying the preceding formula to this identity gives
\[
\begin{split}
 D^j\chi(y)[v_1,\ldots,v_j]
 ={}&-\bigl(D\eta(\chi(y))\bigr)^{-1}\\
 &\quad\cdot\sum_{\substack{\pi\in\Pi_j,\ |\pi|\geq2\\
                            \pi=\{B_1,\ldots,B_a\}}}
 D^a\eta(\chi(y))
 [D_{B_1}\chi(y)[v],\ldots,D_{B_a}\chi(y)[v]].
\end{split}
\]
The term for the partition with a single block is
$D\eta(\chi(y))[D^j\chi(y)[v_1,\ldots,v_j]]$; solve for it
using $D^j\operatorname{id}=0$ for $j\geq2$. In every other
partition each block has at most $j-1$ elements. Thus every
derivative of $\chi$ on the right has order less than $j$.
Substituting its earlier formulas leaves factors
$D^\ell\eta\circ\chi$, with $2\leq\ell\leq j$, and bounded
matrices $(D\eta)^{-1}\circ\chi$. Induction in the formula shows
that the weights $\ell-1$ sum to $j-1$. Indeed, the outer factor contributes weight $a-1$, and the
expressions already obtained for the blocks contribute
$\displaystyle\sum_{\ell=1}^{a}(|B_\ell|-1)=j-a$. Their sum
is $j-1$. The same H\"older and change of variables estimates
bound the inverse derivatives.
For Sobolev $\eta$, justify the calculation by regularizing its
expression in smaller charts. The approximations converge in $C^1$,
retain the lower differential bound, and have local inverses on still
smaller common domains. The preceding formulas give uniform bounds.
Convergence of the inverses in $C^1$ and weak compactness in the
spaces $L^q(\mathbb R^m)$, with $1<q<\infty$, identify their limits with the
distributional derivatives of $\chi$.

The same formulas give continuity. Change of variables bounds
composition uniformly in every $L^q(\mathbb R^m)$ used. For a fixed $L^q(\mathbb R^m)$ factor,
approximate it by a smooth compactly supported function; uniform
convergence of the mappings and the Jacobian bound give convergence
of the composed factor in $L^q(\mathbb R^m)$. Split varying factors into their
difference and a fixed factor. H\"older's inequality then controls
differences of the preceding products. When some coefficients are only
bounded, first approximate an integrable factor by one of compact
support: continuity of the coefficients suffices on that compact set,
and their uniform bound controls the tail. This proves both joint
continuity of composition on sections and continuity of the integrable
terms obtained by applying the inversion formula to the displacement
near the identity.

Near the identity these formulas also control the displacement, which
must be integrable. If $\eta(x)=\exp_x\mathbf X(x)$ and
$\chi=\eta^{-1}$, reversing the radial geodesic gives
\[
 \log_y\chi(y)=-P_{\chi(y)\to y}\mathbf X(\chi(y)),\qquad y\in M.
\]
The transport and composition estimates put the right hand side in
$W^{r,p}(M,TM)$, continuously in $\mathbf X$.
For the product of two mappings near the identity, compare their two
successive exponential maps. In radial frames this is a superposition
in their displacements and the composed section; it vanishes when both
sections are zero. The integral difference formula and mixed bounds
through order $r+1$ give membership of the displacement in $W^{r,p}(M,TM)$ and continuity.
Each term thus contains a section factor; no constant is integrated
over $M$. For finitely many factors, approximate their displacements
by smooth compactly supported sections. The products of the
approximations are smooth centers, and the same estimates express the
original product in a chart at one such center. Thus these finite
products still belong to $\Omega^{p,r}(M,M)$.

We now identify the component. Exponential coordinates centered at
$\eta\in\mathcal D^{p,r}(M)$ were justified at the end of the atlas
proof. If $h=\operatorname{Exp}_\eta\mathbf Y$, then
\[
 \log_y\bigl(h(\eta^{-1}(y))\bigr)
 =\mathbf Y(\eta^{-1}(y)),\qquad y\in M.
\]
The composition estimate just proved controls this section in
$W^{r,p}(M,TM)$. There is therefore a neighborhood $\mathcal O_\eta$
on which $\Theta_\eta(h):=h\circ\eta^{-1}$ belongs to an identity
chart and depends continuously on $h$. The inverse composition
$u\mapsto u\circ\eta$ obeys the same estimates. These local
statements do not presuppose that all of $\mathcal D^{p,r}(M)$ is
a group.

Choose a symmetric path-connected neighborhood $U$ of the identity
where the local estimates hold. Let $H$ be the subgroup of the abstract
group of $C^1$ diffeomorphisms generated by $U$. Successive application
of the composition estimates to finitely many factors ensures
$H\subseteq\mathcal D^{p,r}(M)$ and continuity of those compositions.
Each product of elements of $U$ is joined to the identity by multiplying
paths of its factors; hence $H\subseteq\mathcal D_0^{p,r}(M)$.
For each $\eta\in\mathcal D^{p,r}(M)$, shrink
$\mathcal O_\eta$ so that $\Theta_\eta(\mathcal O_\eta)\subseteq U$.
Then $\mathcal O_\eta\subseteq H\eta$. The intersections of right
cosets $H\eta$ with $\mathcal D^{p,r}(M)$ are open and form a partition.
In particular, $H$ is open and closed in $\mathcal D^{p,r}(M)$.
Connectedness gives $\mathcal D_0^{p,r}(M)\subseteq H$, proving equality.
The local estimates, translated by finite products, prove continuity
of composition and inversion on $H$. This gives the topological group
$\mathcal D_0^{p,r}(M)$ with a countable neighborhood base
at the identity. Pairs $(a,b)$ with $a^{-1}b$ in one of these
neighborhoods form a base for its left uniformity. Continuity of
multiplication allows symmetric neighborhoods to be chosen whose cubes lie
in the preceding neighborhood; Proposition~\ref{prop:metrizacion-cadenas-uniformidad}
gives a compatible metric. These are the properties of
Theorems~6.4 and~6.5 of \cite[\S6]{EichhornMaps1993} used here.

The topology retains global control in the Sobolev norms of the charts.
At fixed regularity, one obtains a topological group; the differentiable
properties recovered by moving to higher levels motivate the ILH tower.

\subsection{Hilbert scales and ILH towers}

The preceding loss of derivatives naturally leads to working with a scale of spaces instead of fixing a single regularity.

\begin{definition}[Hilbert scale]
\label{def:escala-hilbert}\glsadd{escala-hilbert}
\index{Hilbert scale}
Let $s_0\in\mathbb Z$. A \textit{Hilbert scale} is a family
$(\mathscr H^s)_{s\in\mathbb Z,\ s\geq s_0}$ with continuous dense inclusions
$\mathscr H^{s+1}\hookrightarrow \mathscr H^s$. Its smooth core is
$\displaystyle \mathscr H^\infty:=\bigcap_{s\in\mathbb Z,\ s\geq s_0}\mathscr H^s$, equipped with the initial topology of the
inclusions $\mathscr H^\infty\longrightarrow \mathscr H^s$; equivalently,
$\mathscr H^\infty=\displaystyle\varprojlim_{s\in\mathbb Z,\ s\geq s_0}\mathscr H^s$.
\end{definition}

\begin{definition}[Differentiable compatibility with loss]
\label{def:compatibilidad-ilh-perdida}
Fix $\ell\in\mathbb N_0$ and $s\in\mathbb Z$ with $s\geq s_0$.
A family of maps
$F^s\colon\mathcal U^{s+\ell}\longrightarrow\mathcal N^s$ is an
\textit{extension between levels of class $\ell$} of $F^\infty$ if its
restrictions are compatible with the inclusions, each $F^s$ is of class
$C^\ell$, and, for $0\leq j\leq\ell$, its differentials, expressed in compatible charts, satisfy on
neighborhoods with common control data an estimate
\[
 \|D^jF^s(x)(v_1,\ldots,v_j)\|_{\mathscr H^s}
 \leq C\,\mathcal P(\|x\|_{\mathscr H^{s+\ell}})
 \prod_{i=1}^j\|v_i\|_{\mathscr H^{s+\ell}},
\]
where $\mathcal P$ is a polynomial independent of the smooth approximation.
The constant may depend on the charts and their control data. For
diffeomorphisms, in addition to Sobolev bounds, a common bound
on $\|(d\eta)^{-1}\|_\infty$ is fixed; boundedness of the $H^s(\Omega,\mathbb R^m)$ norms of the coordinate expressions of $\eta$ on the domains $\Omega$ of the fixed charts
alone does not prevent the differential from approaching a singular matrix.
We say that $F^\infty$ is \textit{smooth in the ILH tower} if it admits these
extensions for every $\ell$. This definition uses Banach calculus
at each level and explicitly records the derivatives that must be recovered by
ascending the scale.
\end{definition}

\begin{definition}[Compatible ILH tower]
\label{def:variedad-grupo-ilh}\glsadd{variedad-ilh}
\index{ILH manifold}
\index{ILH group}
An \textit{ILH tower} modeled on $(\mathscr H^s)$ is a family of Hilbert manifolds
$(\mathcal M^s)_{s\in\mathbb Z,\ s\geq s_0}$ with smooth dense inclusions
$\mathcal M^{s+1}\hookrightarrow\mathcal M^s$ with the following local compatibility. At every point of the core there is
a base-level chart $\varphi^{s_0}\colon U^{s_0}\to V^{s_0}\subset \mathscr H^{s_0}$
whose restrictions, for every $s\geq s_0$, are charts
\[
 \varphi^s\colon U^s=U^{s_0}\cap\mathcal M^s
       \longrightarrow V^s=V^{s_0}\cap \mathscr H^s.
\]
These charts cover the core and their transitions are compatible at all levels.
Thus we require domains coming from one base-level open set, not merely an
arbitrary intersection of open sets. The core is
$\displaystyle \mathcal M^\infty:=\bigcap_{s\in\mathbb Z,\ s\geq s_0}\mathcal M^s$. An ILH group tower
also has a group operation on $G^\infty$ and extensions of
multiplication and inversion between levels. The regularities of these
extensions and the loss of derivatives are specified in
Definition~\ref{def:compatibilidad-ilh-perdida}; the mere equality
$\displaystyle G^\infty=\bigcap_{s\in\mathbb Z,\ s\geq s_0}G^s$ does not suffice.
\end{definition}

\begin{proposition}
\label{prop:nucleo-escala-frechet}
If the inclusions $\mathscr H^{s+1}\hookrightarrow \mathscr H^s$ are continuous, then $\mathscr H^\infty$, with the norms
$\|\cdot\|_{\mathscr H^s}$, is a Fréchet space. A sequence converges in
$\mathscr H^\infty$ if and only if it converges at every level $\mathscr H^s$.
\end{proposition}

\begin{proof}
Take $Y=\mathscr H^{s_0}$ and $X_n=\mathscr H^{s_0+n-1}$ in
Proposition~\ref{prop: completitud secuencial limite proyectivo}.
The inclusions into $Y$ are continuous by composition of successive
inclusions, and $Y$ is Hausdorff. That proposition gives sequential
completeness of the core in its projective topology. The countable
family of norms separates points;
Theorem~\ref{teo: caracterizacion Frechet por seminormas} shows that the
core is Fréchet. The convergence criterion is
Proposition~\ref{prop: convergencia y Cauchy por seminormas}.
\end{proof}

\begin{proposition}[Smooth core of the Sobolev scale]
\label{prop:nucleo-suave-escala-sobolev}
If $M$ is a smooth compact manifold with or without boundary and
$\mathbf{E}\longrightarrow M$ is a smooth vector bundle, then
$\displaystyle \Gamma(\mathbf{E})=\bigcap_{s\in\mathbb Z,\ s\geq s_0}H^s(M,\mathbf{E})$ as locally convex spaces:
the projective topology of the Sobolev norms agrees with the usual topology
of uniform convergence of all derivatives. If $M$ is noncompact,
we always have $\displaystyle \bigcap_{s\in\mathbb Z,\ s\geq s_0}H^s(M,\mathbf{E})\subseteq\Gamma(\mathbf{E})$, but the
inclusion may be strict because a smooth section need not be
integrable at infinity.
\end{proposition}

\begin{proof}
First suppose $M$ is compact. For fixed $s$, choose an integer
$k\geq\displaystyle\max\{s,0\}$. Inclusion in the scale and equivalence between the
$H^k(M,\mathbf E)$ norm and the covariant norm of order $k$ give constants
$A_{s,k},B_k>0$ such that, for $\mathbf{u}\in\Gamma(\mathbf{E})$,
\[
 \|\mathbf{u}\|_{H^s(M,\mathbf{E})}\leq A_{s,k}\|\mathbf{u}\|_{H^k(M,\mathbf{E})}
 \leq B_k\lambda_{\mathbf{g}}(M)^{1/2}
 \sum_{\ell=0}^k
 \|\nabla^\ell \mathbf{u}\|_{
 C^0(M,T^{(0,\ell)}(TM)\otimes \mathbf{E};\mathbf{g},\mathbf{h}_{\mathbf{E}})}.
\]
Thus each inclusion $\Gamma(\mathbf{E})\longrightarrow H^s(M,\mathbf{E})$ is continuous. In the
opposite direction, given $j$, choose $s>j+m/2$. Sobolev embedding gives
$\|\mathbf{u}\|_{C^j(M,\mathbf{E})}\leq C_{s,j}\|\mathbf{u}\|_{H^s(M,\mathbf{E})}$. Thus the scale norms
and the $C^j$ seminorms generate the same topology. On a noncompact
base, even a nonzero parallel section may have infinite $L^2(M,\mathbf E)$
norm, proving the last assertion.
\end{proof}

\begin{definition}[The term \emph{strong ILH}]
\label{def:grupo-strong-ilh}
\index{strong ILH Lie group@\emph{strong ILH} Lie group}
We reserve the term \textit{strong ILH Lie group} for Omori's
notion~\cite[Chapter~I, \S1.1]{Omori1997}. It requires a Hilbert chain, a
compatible chart at the identity, and conditions $(N,1)$--$(N,7)$ for
extensions of multiplication, translations, their iterated
differentials, inversion, and conjugation by smooth elements. A compatible ILH
tower in the sense of the preceding definitions will not be regarded as
\emph{strong ILH} until these axioms have been verified.
\end{definition}

\begin{remark}[Three forms of uniform control]
\label{obs:tres-nociones-uniformes}
\label{obs:distincion-amann-eichhorn-ilh-inicial}
\index{bounded geometry}
\index{uniform structure!Eichhorn}
\index{ILH}
The notions used in this chapter act on different objects.
Amann's uniform regularity is a coordinate property of the atlas and the
base metric; bounded geometry is its intrinsic formulation through
injectivity radius, curvature, and, when there is boundary, collar data. An
Eichhorn uniformity is imposed on a configuration space and serves
to define global closeness, components, and completions. An ILH tower,
finally, organizes a family of Hilbert completions and losses of
regularity between them.

Theorems~\ref{teo:amann-ur-equivale-geometria-acotada}
and~\ref{teo:amann-ur-equivale-geometria-acotada-frontera} relate the first
two notions, without and with boundary. Neither determines by itself
a uniformity on configurations, and a uniformity does not provide
the charts and estimates between levels required by an ILH tower.
\end{remark}

\begin{proposition}[Operations on the diffeomorphism core]
\label{prop:operaciones-nucleo-difeomorfismos-ilh}
Let $M$ be a smooth compact manifold with or without boundary and let
$s_0\in\mathbb Z$ with $s_0>\frac{m}{2}+1$. Then
\[
 \operatorname{Diff}^{\infty}(M)
 =\bigcap_{s\in\mathbb Z,\ s\geq s_0}\operatorname{Diff}^{s}(M),
\]
and composition and inversion are smooth in the sense of
Definition~\ref{def:compatibilidad-ilh-perdida}. More precisely, for each
$\ell\geq0$,
\[
 \operatorname{Diff}^{s+\ell}(M)\times\operatorname{Diff}^{s}(M)
 \longrightarrow\operatorname{Diff}^{s}(M),
 \qquad(\eta,\varphi)\longmapsto\eta\circ\varphi,
\]
and
\[
 \operatorname{Diff}^{s+\ell}(M)\longrightarrow\operatorname{Diff}^{s}(M),
 \qquad\eta\longmapsto\eta^{-1},
\]
are of class $C^\ell$.
\end{proposition}

\begin{proof}
The set-theoretic equality follows from
Proposition~\ref{prop:nucleo-suave-escala-sobolev} in charts, for both
the map and its inverse. Composition from levels $s+\ell$ and $s$ is
$C^\ell$ by
Lemma~\ref{lem:composicion-sobolev-perdida-autocontenida}. Inversion
between those levels is $C^\ell$ by
Theorem~\ref{teo:ebin-marsden-diff-sobolev}, whose proof includes
continuity at the initial level and the implicit function argument.

The extensions agree on the smooth core because they are the same
pointwise operations of composition and inversion. Their differentials
also agree on restriction to directions of higher regularity: the
difference quotients are the same, and continuous inclusions preserve
their limits. Repeat this argument for iterated derivatives.
Let us also verify the local character of the bounds
required by Definition~\ref{def:compatibilidad-ilh-perdida}. Fix a point,
a level and $\ell$. The derivatives of the extensions through order
$\ell$ are continuous in operator norm. On a sufficiently small
neighborhood their norms are bounded by a common constant, taking
the maximum of the $\ell+1$ bounds. This gives the estimate in that
definition with the constant polynomial $\mathcal P=1$ on this
neighborhood. It may also be reduced to fix a bound for
$\|(d\eta)^{-1}\|_\infty$. The explicit composition estimates and the
inversion recurrence preserve bounds on the controlled families
specified there. No common bound on the whole group is being asserted.
This verifies the required regularity and compatibility.
Passage to the Fréchet atlas also uses the common domains of
Corollary~\ref{cor:torre-ilh-variedad-frechet}; it does not follow merely
from the set-theoretic equality of the intersections.
\end{proof}

This definition does not assert that each $G^s$ is a Hilbert Lie group. In
the diffeomorphism example, $G^s$ is a topological group with
smooth right translations; the smooth operation reappears in the limit
because any finite loss can be compensated by ascending finitely many
levels.

Scales do not remove derivative loss; they record it. An argument at
fixed regularity must close in a single $\mathscr H^s$. An ILH argument may temporarily ascend
to $\mathscr H^{s+j}$ and then descend, provided the estimates
are uniform and compatible. This distinction will be decisive when interpreting
the geometric differential operators of the next chapter as fields
on scales of function spaces.

\begin{example}[A smooth isotopy does not determine the Sobolev component]
\label{ej:traslaciones-no-componente-sobolev}
On $M=\mathbb R$, the Euclidean exponential map is addition.
A Sobolev chart centered at a map $f$ therefore lies in the affine
space $f+W^{r,p}(\mathbb R)$, $1<p<\infty$.
Two such affine spaces are either equal or disjoint, and each is open
and closed in the union of the charts: every coordinate neighborhood
stays in the same class modulo $W^{r,p}(\mathbb R)$. In particular, the identity
component of $\mathcal D^{p,r}$ lies in
$\operatorname{id}+W^{r,p}(\mathbb R)$.

If $a\ne0$, the translation $\tau_a(x)=x+a$ is a smooth diffeomorphism,
satisfies $d\tau_a=1$ and has all higher derivatives equal to zero.
The family $\tau_{ta}$, $0\leq t\leq1$, is a smooth isotopy with the
same uniform bounds, and the smallest singular value of its differential
is one. Nevertheless,
\[
 \tau_a-\operatorname{id}=a\notin L^p(\mathbb R),
\]
so $\tau_a$ is not in the Sobolev identity component at any finite $p$
level. Nor does it belong to the intersection of those components for
$p=2$.

The isotopy requirement must be made precise when comparing with
Proposition~8.4 in Chapter III of \cite{Eichhorn2007}: an isotopy in
the class of smooth maps with uniformly bounded derivatives is not
enough for membership in the identity component at every Sobolev level.
If an isotopy continuous at those levels is required, the preceding
translation family is excluded. The topology of the isotopy and
integrability at infinity cannot be omitted when describing the
noncompact ILH core.
\end{example}

\begin{remark}[Hierarchy of structures]
\label{obs:jerarquia-banach-ilh}
A Hilbert manifold at a fixed regularity level, a compatible tower, and a
\emph{strong ILH} Lie group are successively stronger structures. The identity
$\displaystyle G^\infty=\bigcap_{s\in\mathbb Z,\ s\geq s_0}G^s$ alone provides neither compatible charts nor
estimates for derivative loss, and these charts do not automatically imply
Omori's axioms either. In what follows, constructions at fixed regularity are
treated using Banach calculus, while operations with derivative loss are
organized by means of Definition~\ref{def:compatibilidad-ilh-perdida}.
\end{remark}

\section{Finite-regularity maps and manifolds of smooth maps}
\label{sec:mapeos-clasicos-frechet}

Sobolev exponential charts have versions on $C^k$ spaces and on the
space of all smooth maps. The common construction allows the topologies
to be compared without identifying them. In the next propositions, $M$
is a compact finite-dimensional smooth manifold, with or without boundary,
and $N$ is a finite-dimensional smooth manifold without boundary.
Completeness of an auxiliary metric on $N$ is not required.

\begin{lemma}[Superposition on $C^k$ section spaces]

\label{lem:omega-Ck-frechet}

Let $\mathbf E,\mathbf F\to M$ be smooth finite-rank bundles,
let $\mathcal O\subseteq\mathbf E$ be open, and let
$A\colon\mathcal O\to\mathbf F$ be a smooth base-preserving map.
For $k\in\mathbb N_0$, the set

\[
 \mathcal U^k=\{\mathbf u\in C^k(M,\mathbf E)\mid\mathbf u(M)\subseteq\mathcal O\}
\]

is open, and $\Omega_A\colon\mathcal U^k\to C^k(M,\mathbf F)$,
$\Omega_A(\mathbf u)=A\circ\mathbf u$, is smooth. Its derivatives are

\begin{equation}\label{eq:omega-Ck-derivadas}
 D^q\Omega_A(\mathbf u)(\mathbf v_1,\ldots,\mathbf v_q)(x)
 =D_v^qA(x,\mathbf u(x))
                (\mathbf v_1(x),\ldots,\mathbf v_q(x)).
\end{equation}

The restriction to smooth sections is a smooth Fréchet map
$\mathcal U^\infty\to\Gamma(\mathbf F)$ with the same formula.
For finite $k,r$, the assumption $A\in C^{k+r}$ suffices to obtain
$\Omega_A\in C^r$.

\end{lemma}

\begin{proof}

For a base without boundary, Lemma~\ref{lem:algebra-ck-compacta}
directly proves completeness of the covariant norm and the product
estimates, including those for finite-rank sections.

The $C^k$ norms are defined using finitely many charts and
trivializations on compact sets covering $M$.
Proposition~\ref{prop: GammaKm es Banach}, with $K=M$ and order $k$,
proves completeness. Theorem~\ref{norma ck independiente de metrica y conexion}
gives equivalence with the covariant norms and independence of the
auxiliary data.

Fix $\mathbf u\in\mathcal U^k$ and a metric on the fibers.
For $\varepsilon>0$, the set
\[
 T_\varepsilon(\mathbf u)
 :=\{v\in\mathbf E_x\mid x\in M,
                         |v-\mathbf u(x)|\leq\varepsilon\}
\]
is compact: it is the image of the closed disk bundle of radius
$\varepsilon$ under continuous translation by $\mathbf u$.
There exists $\varepsilon>0$ such that
$T_{2\varepsilon}(\mathbf u)\subset\mathcal O$.
Otherwise, we could choose $v_j\notin\mathcal O$ with
$|v_j-\mathbf u(x_j)|\leq1/j$. A subsequence in
$T_1(\mathbf u)$ would converge to $\mathbf u(x)$,
contradicting closedness of $\mathbf E\setminus\mathcal O$.
The ball $\|\mathbf w-\mathbf u\|_{C^0}<\varepsilon$ and all
segments $\mathbf u+t(\mathbf w-\mathbf u)$, $0\leq t\leq1$,
have their graphs in $T_\varepsilon(\mathbf u)$.
Continuity of the inclusion $C^k\hookrightarrow C^0$ proves
openness. On a smaller $C^k$ ball around $\mathbf u$, norms are
bounded and all graphs remain in this same compact subset of
$\mathcal O$.

Work in one trivialization and
write the map as $A(x,u)$. On spatially differentiating the right-hand
side of \eqref{eq:omega-Ck-derivadas}, each term has the form

\[
 (\partial_x^\beta D_v^{q+a}A)(x,u(x))
 \bigl(\partial^{\gamma_1}u,\ldots,\partial^{\gamma_a}u,
        \partial^{\delta_1}v_1,\ldots,\partial^{\delta_q}v_q\bigr),
\]

multiplied by a combinatorial integer. Here $a\geq0$,
$|\gamma_i|\geq1$, $|\delta_i|\geq0$, and
$\beta+\displaystyle \sum_{i=1}^a\gamma_i+\sum_{i=1}^q\delta_i=\alpha$, with
$|\alpha|\leq k$. This description follows by induction: a new derivative
hits the coefficient, a derivative of $u$, or a derivative of one of the
$v_i$, preserving the multi-index sum. If $u$ stays in a bounded $C^k$
ball and its graph stays in a fixed compact subset of $\mathcal O$,
we obtain

\[
 \|D_v^q A(\,\cdot\,,u)(v_1,\ldots,v_q)\|_{C^k}
 \leq C_{k,q,B}\prod_{i=1}^q\|v_i\|_{C^k}.
\]

The spatial formula also proves continuity of these coefficients
in $u$: derivatives of $u$ converge uniformly, and derivatives of $A$
are uniformly continuous on the chosen compact set. Now apply integral
Taylor in the variable $v$ to $D_v^qA$. If $A$ is smooth, the remainder
after the linear term is bounded in $C^k$ by
$C\|h\|_{C^k}^2\displaystyle \prod_{i=1}^{q}\|v_i\|_{C^k}$, using the same estimates through
order $k+q+2$. Induction proves that
\eqref{eq:omega-Ck-derivadas} are operator-norm Fréchet derivatives,
not merely directional formulas.

If only $A\in C^{k+r}$, use in the last step the first-order
remainder written as an integral of a derivative difference.
Given $u$ and $\varepsilon>0$, uniform continuity of these derivatives
on the chosen compact set provides $\delta>0$ such that the remainder
has $C^k$ norm at most
$\varepsilon\|h\|_{C^k}\displaystyle\prod_{i=1}^{r-1}\|v_i\|_{C^k}$
if $\|h\|_{C^k}<\delta$. The empty product equals one when $r=1$.
This bound proves the derivative of order $r$ and its continuity without
requiring an additional derivative of $A$.
Finally, on smooth sections, each $C^k$ seminorm is controlled by the
corresponding Banach assertion. Difference quotients converge in each
seminorm, and the derivative evaluation formulas are jointly continuous
for each one. This is precisely Bastiani smoothness.

\end{proof}

\begin{lemma}[Approximation by smooth maps]

\label{lem:aproximacion-Ck-mapeos}

For every $k\in\mathbb N_0$, smooth maps $M\to N$ are dense in
$C^k(M,N)$. The approximation may be chosen arbitrarily close in $C^0$.

\end{lemma}

\begin{proof}
Fix $f\in C^k(M,N)$. Choose coordinate domains $U_i$ in $M$,
compact sets $K_i\subseteq U_i$ covering $M$, and target charts
$\psi_i\colon V_i\to\mathbb R^n$ such that
$f(\overline{U_i})$ is a compact subset of $V_i$.
Choose cutoffs $\chi_i$ supported in $U_i$ and equal to one near
$K_i$. Write $1\leq i\leq m$ for this finite cover.

Choose open sets $V'_i$ with
$f(\overline{U_i})\subset V'_i\Subset V_i$ and a compatible
distance $d_N$. Compactness and finiteness give $\delta_0>0$ such
that
\[
 x\in\overline{U_i},\quad d_N(y,f(x))<\delta_0
 \quad\Longrightarrow\quad y\in V'_i
 \qquad(1\leq i\leq m).
\]
Each compact set $\psi_i(\overline{V'_i})$ has a positive
Euclidean margin $b_i>0$ from the complement of $\psi_i(V_i)$.

Start with $f_0=f$. At step $i$, approximate
$u_i:=\psi_i\circ f_{i-1}$ on the support of $\chi_i$ by a
smooth function $a_i$ using convolution. An additional cutoff equal
to one on this support allows convolution on a larger domain.
In a boundary chart, first extend the $C^k$ components across the
half-space. Finite reflection in a smaller strip gives such an extension:
for $t<0$,
use $\sum_{j=1}^{k+1}c_j u_i(x',-jt)$, where
$\sum_{j=1}^{k+1}c_j(-j)^\ell=1$ for $0\leq\ell\leq k$.
The Vandermonde system is invertible, and these identities match
derivatives through order $k$ at $t=0$.
Uniform continuity of derivatives on compact sets proves
$a_i\to u_i$ in $C^k$ on the chosen support.

Define
\[
 f_i=\psi_i^{-1}\bigl((1-\chi_i)u_i+\chi_i a_i\bigr)
 \quad\hbox{on }U_i,\qquad
 f_i=f_{i-1}\quad\hbox{outside }U_i.
\]
Choose the approximation with $\|a_i-u_i\|_{C^0}<b_i/2$ and,
using uniform continuity of $\psi_i^{-1}$ on the compact set formed
by these segments, with
\[
 \sup_{x\in M}d_N(f_i(x),f_{i-1}(x))<\frac{\delta_0}{2m}.
\]
Inductively, $d_N(f_i(x),f(x))<\delta_0/2$.
Thus all subsequent steps are defined in their target charts, and
the preceding coordinate segments stay in $\psi_i(V_i)$.
The expressions glue because $\chi_i$ vanishes near the edge of
$U_i$.

The map $f_i$ is smooth near $K_i$ and preserves smoothness wherever
$f_{i-1}$ was already smooth. Hence $f_m$ is smooth on all of $M$.
Superposition by smooth changes of coordinates is continuous in
$C^k$, by Lemma~\ref{lem:omega-Ck-frechet} on the preceding compact
sets. Given a $C^k$ neighborhood of $f$, describe it by finitely many
coordinate seminorms and positive bounds. At each step, choose the
error in each such seminorm smaller than $1/m$ of its bound.
The triangle inequality places $f_m$ in the prescribed neighborhood
while retaining the preceding $C^0$ bounds.
\end{proof}

\begin{theorem}[A common atlas for $C^k(M,N)$ and $C^\infty(M,N)$]

\label{teo:variedades-mapeos-Ck-frechet}

For $k\in\mathbb N_0$, $C^k(M,N)$ is a smooth Banach manifold.
The space $C^\infty(M,N)$ is a smooth Fréchet manifold.
Both structures are constructed with the same exponential charts,
restricted to sections of the corresponding regularity. In particular,

\[
 T_fC^k(M,N)=C^k(M,f^*TN),\qquad
 T_fC^\infty(M,N)=\Gamma(f^*TN).
\]

The smooth structures are independent of the auxiliary metric used
to construct the charts.

\end{theorem}

\begin{proof}
Fix a metric $\mathbf h$ on $N$ and write
$\mathscr E(y,v)=(y,\exp_yv)$. For every compact $K\subset N$,
there exist an open set $O\supset K$ and $r>0$ such that
\[
 \mathscr E\colon
 D_{2r}(O):=\{(y,v)\in TN\mid y\in O,\ |v|_{\mathbf h}<2r\}
 \longrightarrow\mathscr E(D_{2r}(O))
\]
is a diffeomorphism onto an open subset of $N\times N$.
Indeed, the inverse function theorem at each $(y,0)$ gives an
open set $O_y$ and $\rho_y>0$ on which $\mathscr E$ is a
diffeomorphism on $D_{\rho_y}(O_y)$. Finitely many of these open
sets cover $K$. Take $O$ to be their union and
$2r<\min\rho_y$ over this cover. The differential is invertible
on $D_{2r}(O)$. If two points have the same image, their base points
coincide, and both lie in one of the domains
$D_{\rho_y}(O_y)$ where injectivity is already known.
This proves the assertion and fixes a common radius around $K$.

For each smooth center $f\colon M\to N$, include all choices of
an open set $O_f\supset f(M)$ and a radius $\rho_f>0$ for which
$\mathscr E$ is a diffeomorphism on $D_{\rho_f}(O_f)$.
Define
\begin{align*}
 \mathcal V_f^k&=\{\mathbf X\in C^k(M,f^*TN)\mid
                                  \|\mathbf X\|_{C^0}<\rho_f\},\\
 \mathcal U_f^k&=\{g\in C^k(M,N)\mid
           (f(x),g(x))\in\mathscr E(D_{\rho_f}(O_f))
                         \text{ for every }x\in M\},\\
 \operatorname{Exp}_f(\mathbf X)(x)&=\exp_{f(x)}\mathbf X(x),\\
 \Phi_f(g)(x)&=\log_{f(x)}g(x).
\end{align*}
For $k=\infty$, use $\Gamma(f^*TN)$.
The inclusion into $C^0$ is continuous, so $\mathcal V_f^k$ is
open. For $g\in\mathcal U_f^k$, the graph
$x\mapsto(f(x),g(x))$ is a compact subset of the open set
$\mathscr E(D_{\rho_f}(O_f))$. A sufficiently small $C^0$
neighborhood preserves this inclusion, so $\mathcal U_f^k$ is open
as well. The exponential and logarithm formulas are continuous
inverses in $C^k$, by smooth superposition on finitely many compact
coordinate sets and Lemma~\ref{lem:omega-Ck-frechet}.

To prove that these charts cover the space, let $g\in C^k(M,N)$
and apply the common-radius construction to $K=g(M)$.
The open set $W:=\mathscr E(D_r(O))$ contains the diagonal over
$K$. For a compatible distance $d_N$, there exists $\delta>0$
such that
\[
 z\in K,\quad d_N(y,z)<\delta
 \quad\Longrightarrow\quad (y,z)\in W.
\]
Otherwise there would be $z_j\in K$ and $y_j$ with
$d_N(y_j,z_j)\to0$ and $(y_j,z_j)\notin W$. A subsequence
would converge to $(z,z)\in W$, a contradiction.
By Lemma~\ref{lem:aproximacion-Ck-mapeos}, choose a smooth center
$f$ with $\sup_x d_N(f(x),g(x))<\delta$.
Then $f(M)\subset O$, and
$\mathbf X(x)=\log_{f(x)}g(x)$ is a $C^k$ section with
$\|\mathbf X\|_{C^0}<r$. The chart centered at $f$ with radius
$\rho_f=2r$ contains $g$.
For $C^\infty(M,N)$, use the map itself as center.
All radii are chosen over compact subsets of the target.

On an overlap, the transition is superposition by the smooth
bundle map

\[
 (x,v)\longmapsto
 \bigl(x,\log_{g(x)}(\exp_{f(x)}v)\bigr).
\]

Lemma~\ref{lem:omega-Ck-frechet}, applied around the graph of the
section corresponding to the point of overlap, gives a common domain
and local derivative bounds. It follows that the overlap is open and
the transition and its inverse are smooth. Two distinct maps are separated by a point
evaluation and disjoint open subsets of $N$, proving the Hausdorff
property. Local $C^k$ norms, and all the $C^j$ seminorms in the smooth
case, show that the topology obtained is the usual convergence of
components and derivatives on finitely many coordinate compact sets.

Differentiating a curve $s\mapsto f_s$ in a chart and evaluating at
a point gives the section $x\mapsto\partial_sf_s(x)|_{s=0}$ of $f^*TN$.
Conversely, a section $\mathbf V$ of the stated regularity is represented
by $s\mapsto\operatorname{Exp}_f(s\mathbf V)$ when the center is smooth.
At a $C^k$ center, transport from a smooth-centered chart using the
differential of the exponential. This provides $C^k$ trivializations
of $f^*TN$ and the same tangent description. Their transitions are exactly
the derivatives of chart transitions, and are linear topological
isomorphisms by the superposition lemma. This proves the identifications.
Finally, changing the metric replaces the preceding transition by
$\log^{\widetilde{\mathbf h}}_g\circ\exp^{\mathbf h}_f$, another smooth
superposition near a compact graph. This proves metric independence.

\end{proof}

\begin{proposition}[Evaluation and composition on the smooth mapping space]

\label{prop:operaciones-mapeos-frechet}

Let $M,S$ be smooth compact manifolds, with or without boundary, and let
$N,P$ be smooth finite-dimensional manifolds without boundary.
Evaluation $C^\infty(M,N)\times M\to N$ is smooth.
If $H\colon N\to P$ is smooth, postcomposition $f\mapsto H\circ f$
is smooth on all of $C^k(M,N)$ and all of $C^\infty(M,N)$, without
global bounds on the derivatives of $H$. If $a\colon S\to M$ is smooth
and $S$ is compact, precomposition by fixed $a$ is smooth on these spaces.
If, in addition, $M$ has no boundary, the map

\[
 C^\infty(M,N)\times C^\infty(S,M)\longrightarrow C^\infty(S,N),
 \qquad(f,a)\longmapsto f\circ a,
\]

is smooth. In coordinates, its differential and the differential of
evaluation are

\[
 D\operatorname{Comp}(f,a)(\mathbf U,\mathbf V)
 =\mathbf U\circ a+(df\circ a)\mathbf V,
 \qquad
 D\operatorname{ev}(f,x)(\mathbf U,v)=\mathbf U(x)+df_xv.
\]

\end{proposition}

\begin{proof}

For postcomposition, each map $f$ has compact image. A neighborhood
of that compact set with compact closure allows the superposition lemma
to be applied to $H$. Bounds are therefore needed locally near each $f$,
not uniformly on the whole mapping space. For fixed precomposition,
exponential charts represent the operation by the linear pullback of
sections $\mathbf X\mapsto\mathbf X\circ a$. Each $C^j$ seminorm is
controlled by the chain rule and finitely many derivative bounds for $a$
on $S$. Thus the pullback is continuous linear, and its derivatives of
order greater than one vanish.

For joint composition, the assumption $\partial M=\varnothing$
allows Theorem~\ref{teo:variedades-mapeos-Ck-frechet} to be applied
to $C^\infty(S,M)$. Work on finitely many charts, shrunk so
that their images remain in the target charts under $C^0$ perturbations.
Write the operation there as $u\circ v$. A derivative in the outer
variable replaces $u$ by its variation; a derivative in the inner
variable introduces one additional spatial derivative of $u$.
After $q$ functional derivatives and $j$ spatial derivatives, only
derivatives of $u$ through order $j+q$ and derivatives of $v$ and the
variations through order $j$ occur. Leibniz rules give finite products,
continuous from $C^{j+q}\times C^j$ to $C^j$.

To check that these expressions are differentials, use

\[
 u(v+w)-u(v)-Du(v)w
 =\int_0^1\bigl(Du(v+tw)-Du(v)\bigr)w\,dt.
\]

Let $C_j$ be a $C^j$ multiplication constant for the chosen charts.
Given $u,v$, and $\varepsilon>0$, uniform continuity of derivatives of
$u$ through order $j+1$ on the relevant compact sets gives $\delta>0$ such that
\[
 \sup_{0\leq t\leq1}\|Du(v+tw)-Du(v)\|_{C^j}
 \leq\varepsilon/C_j,\qquad \|w\|_{C^j}<\delta.
\]
Reduce the radius to keep the segments in the charts. The $C^j$ norm
of the remainder is bounded by $\varepsilon\|w\|_{C^j}$. The mixed term containing the outer variation
and $w$ is bounded by the product of their $C^{j+1}$ and $C^j$ norms.
Repeating this argument for the derivative expressions proves the
formulas at every order. Each output seminorm requires only finitely many
input seminorms, so derivative evaluations are jointly continuous in
Fréchet topology. Evaluation is proved in the same way, taking a point
as the inner variable; compactness gives the same local derivative controls.

\end{proof}

\begin{theorem}[The group of smooth diffeomorphisms]

\label{teo:diff-suave-frechet}

If $M$ is compact without boundary, $\operatorname{Diff}^\infty(M)$
is open in $C^\infty(M,M)$ and, with composition and inversion, is a
Fréchet Lie group. For a variation $\mathbf V$ of $\eta$, one has

\[
 D\operatorname{Inv}(\eta)\mathbf V
 =-d\eta^{-1}\bigl(\mathbf V\circ\eta^{-1}\bigr).
\]

\end{theorem}

\begin{proof}
Estimate~\eqref{eq:inyectividad-uniforme-perturbacion-compacta}
and the control of $P_\lambda$ in
Section~\ref{sec:grupos-difeomorfismos-ilh} give a $C^1$
neighborhood of each diffeomorphism consisting of injective maps
with invertible differential. The open-and-closed image argument,
with one point chosen in each component, proves surjectivity.
The inverse is smooth by the finite-dimensional inverse function
theorem. Continuity of $C^\infty\to C^1$ gives the asserted
openness. Multiplication is smooth by the preceding proposition.

First prove continuity of inversion. If $\eta_n\to\eta$ in
$C^\infty$, their inverses converge uniformly: for $y\in M$, the
distance between $\eta(\eta_n^{-1}(y))$ and $y$ is bounded by the
uniform distance between $\eta$ and $\eta_n$, and uniform continuity
of $\eta^{-1}$ applies. In charts,

\[
 D\eta_n^{-1}=(D\eta_n)^{-1}\circ\eta_n^{-1}.
\]

Matrix inversion and uniform convergence give convergence in $C^1$.
Differentiating this identity $j$ times expresses derivatives of order
$j+1$ of the inverse through derivatives of $\eta_n$ up to order $j+1$,
inverses of $D\eta_n$, and lower derivatives of the inverse.
Induction gives convergence in every $C^j$ norm.

Fix an exponential chart with smooth center $\eta_*$ and write
$\eta=\operatorname{Exp}_{\eta_*}(\mathbf X)$.
For a direction $\mathbf Y$ in this chart, set
\[
 \eta_t=\operatorname{Exp}_{\eta_*}(\mathbf X+t\mathbf Y),
 \qquad
 \mathbf V=D\operatorname{Exp}_{\eta_*}(\mathbf X)\mathbf Y.
\]
This family is smooth in $(t,x)$ and stays in the neighborhood of
diffeomorphisms for sufficiently small $|t|$.
The finite-dimensional map $F(t,x)=(t,\eta_t(x))$ has differential
\[
 dF_{(0,x)}(a,w)=(a,a\mathbf V(x)+d\eta_xw)
\]
at $(0,x)$, and this differential is invertible.
For each $y\in M$, the inverse function theorem gives an interval
$I_y$ about zero and an open set $O_y$ about $y$ on which the
branch $(t,z)\mapsto\eta_t^{-1}(z)$ is smooth.
Choose finitely many $O_y$ covering $M$ and $\varepsilon>0$ with
$(-\varepsilon,\varepsilon)\subset I_y$ for every selected
interval. The branches agree on overlaps by injectivity of
$\eta_t$. Thus the inverse is smooth on the common product
$(-\varepsilon,\varepsilon)\times M$.

In the output chart centered at $\eta^{-1}$, write
$\mathbf W(t)(y)=\log_{\eta^{-1}(y)}\eta_t^{-1}(y)$.
Reduce $\varepsilon$ so that this formula is defined on all of
$M$, using the continuity of inversion already proved.
Smoothness on the product and compactness of $M$ give, in every
$C^j$ seminorm,
\[
 \frac{\mathbf W(t)-\mathbf W(0)}{t}-\partial_t\mathbf W(0)
 =\int_0^1
   \bigl(\partial_t\mathbf W(st)-\partial_t\mathbf W(0)\bigr)\,ds
 \longrightarrow0.
\]
This verifies the directional derivative along the line
$\mathbf X+t\mathbf Y$ in a fixed chart.
Differentiating $\eta_t\circ\eta_t^{-1}=\operatorname{id}$ and
returning to the intrinsic identification of tangents yields
$D\operatorname{Inv}(\eta)\mathbf V
=-d\eta^{-1}(\mathbf V\circ\eta^{-1})$.

The formula is jointly continuous by the continuity of inversion
just proved and the properties of composition, multiplication, and spatial
differentiation on $C^\infty$. Thus inversion is $C^1_{\mathrm B}$.
If it is already $C^r_{\mathrm B}$, the same formula and smoothness of
those operations make its derivative evaluation $C^r_{\mathrm B}$;
it is therefore $C^{r+1}_{\mathrm B}$. Induction proves smoothness.

\end{proof}

When $M$ has boundary, openness in all maps $M\to DM$ is not
asserted. Use the charts already constructed with a product metric in
a collar of the double and sections whose boundary restrictions are
tangent to the boundary. This space of smooth sections is closed and
complemented: normal trace has a continuous right inverse obtained by
collar extension followed by a cutoff. The estimates hold in every
$C^j$ seminorm. Small $C^1$ perturbations preserve the interior side
by positivity of the normal derivative, as proved for the charts of
$\operatorname{Diff}^s(M)$. Composition is checked directly in these
charts. In boundary coordinates, work in convex half-balls: the segment
between two nearby values of the inner map remains in the half-space.
The integral Taylor formula and the estimates in the proof of
Proposition~\ref{prop:operaciones-mapeos-frechet} apply there to
derivatives continuous up to the boundary. Each output seminorm is
still controlled by finitely many input seminorms. Variation curves
arising from the charts preserve the boundary and the interior side.
For inversion, the inverse function theorem with parameters applies in
charts of the double; positivity of the normal derivative ensures that
the local inverses restrict to the same side. The preceding induction
on seminorms remains valid. This gives the Fréchet Lie group structure
for smooth diffeomorphisms carrying the boundary to itself.

In particular, writing this group as
$\operatorname{Diff}^\infty_{\partial}(M)$, if $N$ has no boundary,
the action
\[
 C^\infty(M,N)\times\operatorname{Diff}^\infty_{\partial}(M)
 \longrightarrow C^\infty(M,N),\qquad(f,\eta)\longmapsto f\circ\eta,
\]
is smooth by the same half-ball estimates. Its differential is
$\mathbf U\circ\eta+(df\circ\eta)\mathbf V$, where
$\mathbf V(x)\in T_{\eta(x)}\partial M$ for $x\in\partial M$.

\subsection{What is obtained on the core of a scale}
\label{subsec:nucleos-escalas-frechet}

\begin{proposition}[A levelwise smoothness criterion]

\label{prop:criterio-suavidad-niveles-frechet}

Let $(E_j)_{j\in\mathbb N_0}$ and $(F_j)_{j\in\mathbb N_0}$ be chains
of Banach spaces with continuous linear inclusions
$E_{j+1}\hookrightarrow E_j$ and $F_{j+1}\hookrightarrow F_j$, and let $E_\infty=\displaystyle \bigcap_{j=0}^{\infty} E_j$ and
$F_\infty=\displaystyle \bigcap_{j=0}^{\infty} F_j$ carry their projective topologies.
These cores are Fréchet spaces. Let $U\subseteq E_\infty$ be open and
let $A\colon U\to F_\infty$. Suppose that, for every $x\in U$, output
level $j$, and order $r\in\mathbb N_0$, there are a level $m$, an open
set $V\subseteq E_m$ containing $x$, and a $C^r$ map
$A_{j,r}\colon V\to F_j$ whose restriction to $V\cap E_\infty$
agrees with $A$, with $V\cap E_\infty\subseteq U$.
Then $A$ is Fréchet smooth.

\end{proposition}

\begin{proof}

Completeness is proved as in
Proposition~\ref{prop:nucleo-escala-frechet}: a Cauchy sequence converges
in every Banach space, and continuity of the inclusions makes the limits
compatible. Density is not used in this argument. Taking finite maxima
of norms gives increasing seminorms when the original norms are not
increasing.

For $v\in E_\infty$, the difference quotients of $A$ converge in
$F_j$ to the differential of an extension $A_{j,1}$. Comparing two
levels, the continuous inclusion carries the higher-level limit to the
lower-level limit; uniqueness of limits proves compatibility.
The limits therefore form an element of $F_\infty$.
Joint continuity of derivative evaluation is checked in each output
seminorm using the corresponding extension. The same argument for
difference quotients of iterated derivatives proves, by induction on
$r$, existence and joint continuity of all derivatives. The level $m$
may increase with both $j$ and $r$. No tame estimates or additional
strong ILH group axioms have been required.

\end{proof}

\begin{corollary}[Common charts and the Fréchet structure of the core]

\label{cor:torre-ilh-variedad-frechet}

A tower with common charts in the precise sense of
Definition~\ref{def:variedad-grupo-ilh} has a Fréchet manifold as its
core. Its loss-compatible operations are Fréchet smooth when they satisfy
the preceding criterion. In particular, the exponential charts of
$C^\infty(M,N)$ are the common restrictions of the $C^k$ and sufficiently
high-order Sobolev charts.

\end{corollary}

\begin{proof}

A common chart has domain $U^0$ at a base level and represents
$U^j=U^0\cap\mathcal M^j$ on $V^j=V^0\cap E_j$.
Intersecting gives a bijection
$U^\infty\to V^0\cap E_\infty$. The latter set is open in the projective
topology because $V^0$ is open at the base level. Transitions and their
inverses satisfy the levelwise criterion and are therefore smooth.
The topology is Hausdorff through inclusion in the base level, and the
charts cover the core by hypothesis.

For mapping spaces, the common domain is
$\|\mathbf X\|_{C^0}<\rho_f$. Sobolev embeddings make it open at all
levels used. Proposition~\ref{prop:nucleo-suave-escala-sobolev} identifies
the section-space core and its topology with $\Gamma(f^*TN)$.
This proves compatibility of structures, not merely equality of their
underlying sets.

\end{proof}

The common-domain condition is essential to this argument.
A countable intersection of open sets need not be open:
$\displaystyle \bigcap_{j\geq1}(-1/j,1/j)=\{0\}$ in $\mathbb R$.
Thus a set-theoretic equality
$\mathcal M^\infty=\displaystyle \bigcap_{j\in\mathbb Z,\ j\geq s_0}\mathcal M^j$ is not by itself a
construction of a Fréchet manifold.

\begin{remark}[A noncompact source requires a different topology choice]

\label{obs:compacto-abierto-vs-eichhorn}

The compact-source case above and Eichhorn's treatment
\cite{Eichhorn2007,EichhornMaps1993} serve different purposes.
On $\mathbb R$, a nonzero constant function is smooth but does not
belong to $L^2(\mathbb R)$. Moreover, for nonzero
$\chi\in C_c^\infty(\mathbb R)$, the translates $\chi(\,\cdot\,-n)$
converge to zero with all derivatives on compact sets, while their global
$H^s(\mathbb R)$ norms remain constant. Thus the compact-open smooth topology and
the topology on the core of all global Sobolev spaces are not interchangeable.

The exponential chart does not transfer automatically either.
For the constant map $f\colon\mathbb R\to S^1$, $f(x)=1$, the normal
domain $\{u\in C^\infty(\mathbb R)\mid |u(x)|<\pi\text{ for every }x\}$
is not open in the compact-open smooth topology: every zero neighborhood
allows a bump of height $2\pi$ outside the compact sets it controls.
The uniform domains and Sobolev components constructed earlier address a
different problem, namely control at infinity.

We therefore retain both constructions: for compact sources, the
common $C^k$, Sobolev, and Fréchet atlas; for noncompact sources, the
bounded-geometry hypotheses, uniformly controlled centers, and specified
functional components. The first case does not replace the second, and
the greater generality of the second does not remove the need to choose
an appropriate energy level for a variational problem.

\end{remark}

\chapter{Geometric applications in infinite dimensions}
\label{cap:aplicaciones-geometricas-inf}

The constructions of the preceding chapter take on their clearest meaning when an equation can be recognized as a geometric condition on a configuration space. In this chapter we apply this viewpoint to harmonic maps, the Yang--Mills equations, and the Euler equation for an incompressible fluid. Although the problems differ, in each case the unknown is no longer viewed as a collection of components, but as a point of an infinite-dimensional manifold.

The procedure will be the same in each application. First we identify the configuration space and its tangent space; we then introduce the natural functional or metric and compute its variation in an arbitrary direction. Integration by parts then reveals the corresponding differential equation. Symmetries are also part of the problem: invariance under reparametrizations or gauge transformations produces degenerate directions and requires us to distinguish the geometric equation from a particular choice of coordinates or gauge.

The \(L^2\) metrics used in these applications are generally weak with respect to the chosen Sobolev topologies. Thus expressions such as ``gradient of the energy'' will initially be understood in the formal \(L^2\) sense: the resulting operator usually loses derivatives and does not by itself define a vector field on the configuration space at that regularity. When an evolution equation or a smooth spray is needed, we shall specify the regularity scale or elliptic mechanism that recovers the required derivatives. The geometric interpretation of these equations therefore depends on estimates controlling the loss and recovery of derivatives.

\begin{semblanzaHistorica}{Three problems, three geometries}
The theory of harmonic maps developed by extending the Dirichlet energy to maps between manifolds; Yang--Mills theory made the curvature of a connection a central functional of gauge theory; Arnold and, in a precise analytic framework, Ebin and Marsden interpreted the Euler equations as geodesics on a diffeomorphism group. These three examples show that the geometry of the configuration space is not merely decorative: it determines the variables, the symmetries, and the regularity at which the equation is well posed~\cite{EbinMarsden1970}.
\end{semblanzaHistorica}

\section{Harmonic maps}
\label{sec:aplicacion-mapas-armonicos}
\index{harmonic map}
\index{Dirichlet energy}
\index{Jacobi operator}

The first example extends the Dirichlet functional from the calculus of variations
to maps between manifolds. Let \((M,\mathbf{g})\) be a compact
Riemannian manifold of dimension $m$ with or without boundary and \((N,\mathbf{h})\) a
finite-dimensional Riemannian manifold without boundary. For a smooth map
\(f\colon M\longrightarrow N\), the differential is a section of
\(T^*M\otimes f^*TN\), and its \textit{Dirichlet energy} is
\begin{equation}
\label{eq:energia-dirichlet-mapas}
 \mathcal E(f):=\frac{1}{2}\int_M |df|_{\mathbf{g},\mathbf{h}}^{2}\,d\lambda_{\mathbf{g}}.
\end{equation}

The Levi--Civita connection of \(\mathbf{h}\) induces a connection, still denoted by
\(\nabla\), on \(f^*TN\) and on \(T^*M\otimes f^*TN\). The
\textit{tension field} of \(f\) is
\begin{equation}
\label{eq:tension-mapa}
 \boldsymbol{\tau}(f):=\operatorname{tr}_{\mathbf{g}}(\nabla df)
 =\sum_{i=1}^m\bigl(
   \nabla_{\mathbf{e}_i}df(\mathbf{e}_i)
   -df(\nabla^M_{\mathbf{e}_i}\mathbf{e}_i)
 \bigr),
\end{equation}
where \((\mathbf{e}_i)_{i=1}^m\) is any local orthonormal frame. This expression is
tensorial and therefore independent of the chosen frame.

\begin{proposition}[First variation of the energy]
\label{prop:primera-variacion-mapas-armonicos}
Let \((M,\mathbf{g})\) be a compact Riemannian manifold of dimension $m$, with or without
boundary, and \((N,\mathbf{h})\) a Riemannian manifold without boundary. Let \(f_t\) be a
smooth variation of \(f=f_0\), and let
\[
 \mathbf{V}:=\left.\partial_t f_t\right|_{t=0}\in\Gamma(f^*TN).
\]
Then
\begin{equation}
\label{eq:primera-variacion-dirichlet}
 \left.\frac{d}{dt}\right|_{t=0}\mathcal E(f_t)
 =-\int_M\langle\boldsymbol{\tau}(f),\mathbf{V}\rangle_{\mathbf{h}}\,d\lambda_{\mathbf{g}}
 +\int_{\partial M}\langle df(\boldsymbol{\nu}),\mathbf{V}\rangle_{\mathbf{h}}\,d\sigma_{\mathbf{g}},
\end{equation}
where \(\boldsymbol{\nu}\) is the outward unit normal when
$\partial M\neq\varnothing$.
\end{proposition}

\begin{proof}
Consider the map
\[
 F\colon(-\varepsilon,\varepsilon)\times M\longrightarrow N,
 \qquad
 F(t,x):=f_t(x).
\]
Since the Levi--Civita connection of \(N\) is torsion-free and the vector fields
\(\boldsymbol{\partial}_t\) and \(\mathbf{e}_i\) commute,
\[
 \nabla_{\boldsymbol{\partial}_t}df_t(\mathbf{e}_i)
 =
 \nabla_{\mathbf{e}_i}\partial_t f_t.
\]
Evaluating at \(t=0\) gives
\[
 \left.\frac{d}{dt}\right|_{t=0}\frac{1}{2}|df_t|_{\mathbf{g},\mathbf{h}}^2
 =
 \sum_{i=1}^m
 \langle\nabla_{\mathbf{e}_i}\mathbf{V},df(\mathbf{e}_i)\rangle_{\mathbf{h}}.
\]

Define the \(1\)-form
\[
 \boldsymbol{\alpha}(\mathbf{X}):=\langle \mathbf{V},df(\mathbf{X})\rangle_{\mathbf{h}}.
\]
The Leibniz rule and the definition of the tension field give
\[
 \operatorname{div}_{\mathbf{g}}\boldsymbol{\alpha}^\sharp
 =
 \sum_{i=1}^m
 \langle\nabla_{\mathbf{e}_i}\mathbf{V},df(\mathbf{e}_i)\rangle_{\mathbf{h}}
 +
 \langle \mathbf{V},\boldsymbol{\tau}(f)\rangle_{\mathbf{h}}.
\]
Consequently,
\[
 \left.\frac{d}{dt}\right|_{t=0}\mathcal E(f_t)
 =
 \int_M\operatorname{div}_{\mathbf{g}}\boldsymbol{\alpha}^\sharp\,d\lambda_{\mathbf{g}}
 -
 \int_M\langle \mathbf{V},\boldsymbol{\tau}(f)\rangle_{\mathbf{h}}\,d\lambda_{\mathbf{g}}.
\]
Stokes' theorem transforms the first term into
\[
 \int_{\partial M}\boldsymbol{\alpha}(\boldsymbol{\nu})\,d\sigma_{\mathbf{g}}
 =
 \int_{\partial M}\langle \mathbf{V},df(\boldsymbol{\nu})\rangle_{\mathbf{h}}\,d\sigma_{\mathbf{g}}.
\]
\end{proof}

If $\partial M=\varnothing$, or if the trace of $f$ is fixed and
only variations with $\operatorname{Tr}\mathbf V=0$ are allowed, the boundary
term vanishes. In this variational problem, the formal
$L^2(M,f^*TN)$ representative of the differential is
\begin{equation}
\label{eq:gradiente-l2-dirichlet}
 \operatorname{grad}_{L^2(M,f^*TN)}\mathcal E(f)=-\boldsymbol{\tau}(f).
\end{equation}
Consequently, for these variations the critical points are exactly
the maps satisfying $\boldsymbol\tau(f)=0$, called
\textit{harmonic maps}. Variations compactly supported in the
interior establish the necessity of the equation, and the first variation establishes
its sufficiency. If the trace is allowed to vary freely, one must also impose
$df(\boldsymbol\nu)=0$ on the boundary: once the tension field vanishes,
extending any smooth vector field from the boundary makes the
trace of $\mathbf V$ arbitrary in the remaining term. Thus, for
free variations, the boundary term cannot be discarded or generally represented by
an interior $L^2(M,f^*TN)$ gradient.
This variational formulation goes back to
\cite{Eells1966}; the preceding calculation also fixes the sign convention
we shall use.

\begin{proposition}[Linearization of the tension field and second variation]
\label{prop:linealizacion-tension-segunda-variacion}
Let $(M,\mathbf{g})$ be a compact Riemannian manifold of dimension $m$ without boundary,
$(N,\mathbf{h})$ a Riemannian manifold without boundary, and
$f\colon M\longrightarrow N$ a smooth map.
Define
\[
 \nabla^*\nabla:=-\operatorname{tr}_{\mathbf{g}}(\nabla^2)
\]
on sections of \(f^*TN\). The derivative of the tension field is understood
covariantly: transport $\boldsymbol\tau(f_t)$ to $f^*TN$ along
the variation and differentiate at $t=0$. We denote this derivative by
$D\boldsymbol\tau(f)[\mathbf V]$. For every variation with variational vector field
\(\mathbf{V}\),
\begin{equation}
\label{eq:linealizacion-tension-mapa}
 D\boldsymbol{\tau}(f)[\mathbf{V}]
 =
 -\nabla^*\nabla \mathbf{V}
 +\sum_{i=1}^{m}
 \operatorname{Rm}_{\mathbf{h}}\bigl(\mathbf{V},df(\mathbf{e}_i)\bigr)df(\mathbf{e}_i).
\end{equation}
If \(f\) is harmonic, the linearization of \(-\boldsymbol{\tau}\) is the Jacobi operator
\begin{equation}
\label{eq:jacobi-mapa-armonico}
 J_f\mathbf{V}
 =
 \nabla^*\nabla \mathbf{V}
 -\sum_{i=1}^m
 \operatorname{Rm}_{\mathbf{h}}\bigl(\mathbf{V},df(\mathbf{e}_i)\bigr)df(\mathbf{e}_i).
\end{equation}
The Hessian of the energy is given by
\begin{equation}
\label{eq:hessiana-bilineal-mapa-armonico}
 D^2\mathcal E(f)[\mathbf{V},\mathbf{W}]
 =\int_M\biggl(
\langle\nabla \mathbf{V},\nabla \mathbf{W}\rangle_{\mathbf{g},f^*\mathbf{h}}-
 \sum_{i=1}^m
 \left\langle
 \operatorname{Rm}_{\mathbf{h}}\bigl(\mathbf{V},df(\mathbf{e}_i)\bigr)df(\mathbf{e}_i),\mathbf{W}
 \right\rangle_{f^*\mathbf{h}}
 \biggr)\,d\lambda_{\mathbf{g}}.
\end{equation}
In particular,
\begin{equation}
\label{eq:hessiana-mapa-armonico}
 D^2\mathcal E(f)[\mathbf{V},\mathbf{V}]
 =
 \int_M\left(
 |\nabla \mathbf{V}|_{\mathbf{g},f^*\mathbf{h}}^2
 -
 \sum_{i=1}^m
 \left\langle
 \operatorname{Rm}_{\mathbf{h}}\bigl(\mathbf{V},df(\mathbf{e}_i)\bigr)df(\mathbf{e}_i),\mathbf{V}
 \right\rangle_{f^*\mathbf{h}}
 \right)\,d\lambda_{\mathbf{g}}.
\end{equation}
\end{proposition}

\begin{proof}
Let \(F(t,x)=f_t(x)\), and extend the frame \((\mathbf{e}_i)\) locally without
dependence on \(t\). The absence of torsion implies
\[
 \nabla_{\boldsymbol{\partial}_t}df_t(\mathbf{e}_i)=\nabla_{\mathbf{e}_i}\mathbf{V}.
\]
Commuting the covariant derivatives once more gives
\[
 (\nabla_{\boldsymbol{\partial}_t}\nabla_{\mathbf{e}_i}
  -\nabla_{\mathbf{e}_i}\nabla_{\boldsymbol{\partial}_t})df_t(\mathbf{e}_i)
 =
 \operatorname{Rm}_{\mathbf{h}}\bigl(dF(\boldsymbol{\partial}_t),dF(\mathbf{e}_i)\bigr)df_t(\mathbf{e}_i),
\]
so that, at \(t=0\),
\[
 \nabla_{\boldsymbol{\partial}_t}\nabla_{\mathbf{e}_i}df_t(\mathbf{e}_i)
 =
 \nabla_{\mathbf{e}_i}\nabla_{\mathbf{e}_i}\mathbf{V}
 +
 \operatorname{Rm}_{\mathbf{h}}\bigl(\mathbf{V},df(\mathbf{e}_i)\bigr)df(\mathbf{e}_i).
\]
At a given point, choose the frame so that
\(\nabla^M_{\mathbf{e}_i}\mathbf{e}_i=0\). Taking the trace gives
\eqref{eq:linealizacion-tension-mapa}; since both sides are tensorial, the
identity holds at every point and for any frame.

Now suppose that \(\boldsymbol{\tau}(f)=0\). Differentiating the first variation in a
second direction produces a term from the variation of the second vector field;
this term is multiplied by \(\boldsymbol{\tau}(f)\) and vanishes. Thus
\[
 D^2\mathcal E(f)[\mathbf{V},\mathbf{W}]
 =
 \int_M\langle J_f\mathbf{V},\mathbf{W}\rangle_{f^*\mathbf{h}}\,d\lambda_{\mathbf{g}}.
\]
Integrating \(\nabla^*\nabla\) by parts gives
\eqref{eq:hessiana-bilineal-mapa-armonico}. The symmetries of the
curvature tensor show that this form is symmetric. Finally, the principal
symbol of \(J_f\) is
\[
 \boldsymbol\sigma_2(J_f)(x,\xi)
 =-|\xi|_{\mathbf{g}}^2\operatorname{id}_{T_{f(x)}N},
 \qquad x\in M,\quad \xi\in T_x^*M\setminus\{0\},
\]
with the symbol convention from the chapter on differential operators,
which replaces $\partial_j$ by $\xi_j$ without introducing the factor $i$.
The sign comes from the principal term $-g^{ij}\partial_i\partial_j$.
Thus $J_f$ is elliptic; integration by parts and the curvature symmetries
separately establish its formal self-adjointness.
\end{proof}

We turn to the formulation in Sobolev spaces. Fix an integer
\[
 s>\frac{m}{2}+1.
\]
Exponential charts make \(H^s(M,N)\) into a Hilbert manifold.
More precisely, Theorem~\ref{teo:cartas-exponenciales-mapeos-sobolev}
constructs the isomorphism
\[
 \mathcal I_f\colon
 H^s(M,f^*TN)
 \longrightarrow
 T_fH^s(M,N),
 \qquad
 \mathcal I_f(\mathbf{V})
 :=
 \left.\frac{d}{dt}\right|_{t=0}\operatorname{Exp}_f(t\mathbf{V}).
\]
In the chart \(\operatorname{Exp}_f\), this map and its inverse are represented
by the identity on the model space. In particular,
\[
 T_fH^s(M,N)\cong H^s(M,f^*TN).
\]

The tension field contains two derivatives of \(f\), so its natural codomain
has two fewer orders of regularity. To describe it, consider
\[
 \boldsymbol{\mathcal T}^{s-2}
 :=
 \coprod_{\widetilde f\in H^s(M,N)}
 \left(
 \{\widetilde f\}
 \times
 H^{s-2}(M,\widetilde f^*TN)
 \right),
 \qquad
 \pi_{\boldsymbol{\mathcal T}^{s-2}}(\widetilde f,\mathbf{V}):=\widetilde f.
\]
Let \(\mathcal U_f\) be the domain of an exponential chart centered at \(f\). If
\(\widetilde f=\operatorname{Exp}_f(u)\in\mathcal U_f\), parallel
transport along
\[
 t\longmapsto\exp_{f(x)}^N(tu(x))
\]
defines a bundle isometry
\[
 P_{f,\widetilde f}\colon
 \widetilde f^*TN\longrightarrow f^*TN.
\]
The corresponding trivialization is
\[
 \Theta_f\colon
 \pi_{\boldsymbol{\mathcal T}^{s-2}}^{-1}(\mathcal U_f)
 \longrightarrow
 \mathcal U_f\times H^{s-2}(M,f^*TN),
 \qquad
 \Theta_f(\widetilde f,\mathbf{V})
 :=
 \bigl(\widetilde f,P_{f,\widetilde f}\mathbf{V}\bigr).
\]

These trivializations equip \(\boldsymbol{\mathcal{T}}^{s-2}\) with the structure of a smooth
Banach bundle. The analytic point is continuity of multiplication
between the coefficients of the transition maps and the sections.
After multiplying by fixed cutoff functions for a trivialization and extending to
$\mathbb R^m$, write
\[
a\colon\mathbb R^m\to\operatorname{End}(\mathbb R^n),
\qquad
W\colon\mathbb R^m\to\mathbb R^n.
\]
In these representations, we have
\begin{equation}
\label{eq:modulo-sobolev-tension}
 \|aW\|_{H^{s-2}(\mathbb R^m,\mathbb R^n)}
 \leq
 C_s\|a\|_{H^s(\mathbb R^m,\operatorname{End}(\mathbb R^n))}
 \|W\|_{H^{s-2}(\mathbb R^m,\mathbb R^n)},
\end{equation}
where \(a\) is matrix-valued and \(W\) locally represents a section. The
constant depends on \(M\), \(s\), and the trivializations, but not on
\(a\) or \(W\).

To verify the estimate, apply the Leibniz rule to each derivative of
order at most \(s-2\). If $\alpha,\beta\in\mathbb N_0^m$,
$|\alpha|\leq s-2$ and $\beta_i\leq\alpha_i$ for
$i\in\{1,\ldots,m\}$, the regularities of the two factors satisfy
\[
 (s-|\beta|)
 +
 (s-2-|\alpha-\beta|)
 =
 2s-2-|\alpha|
 \geq
 s
 >
 \frac{m}{2}.
\]
Set $a_1=s-|\beta|$ and $a_2=s-2-|\alpha-\beta|$. Both are
nonnegative and $a_1+a_2>m/2$. We may choose
$p_1,p_2\in[2,\infty]$ with $1/p_1+1/p_2=1/2$ so that the
embeddings $H^{a_i}(\mathbb R^m)\hookrightarrow L^{p_i}(\mathbb R^m)$ are continuous. If one of the
$a_i$ is zero, take $p_i=2$ and the other exponent to be infinity, which is allowed
because the other order exceeds $m/2$. If both are positive, the
strict inequality allows the exponents to be chosen away from any
critical endpoint that would require $L^\infty(\mathbb R^m)$. Hölder's inequality then estimates each
product in $L^2(\mathbb R^m)$. Summing over the multi-indices and
the finitely many charts yields \eqref{eq:modulo-sobolev-tension}.

On an overlap of charts centered at \(f_0\) and \(f_1\),
\[
 \Theta_{f_1}\circ\Theta_{f_0}^{-1}(\widetilde f,W)
 =
 \left(
 \widetilde f,
 P_{f_1,\widetilde f}
 P_{f_0,\widetilde f}^{-1}W
 \right).
\]
The coefficient
\[
 \widetilde f
 \longmapsto
 P_{f_1,\widetilde f}P_{f_0,\widetilde f}^{-1}
\]
depends smoothly on $\widetilde f$, with values in
$H^s(M,\operatorname{Hom}(f_0^*TN,f_1^*TN))$, by
Lemma~\ref{lem:superposicion-sobolev-geometrica}. Estimate
\eqref{eq:modulo-sobolev-tension} proves smoothness of the transition
map and its inverse.

To see which terms must be estimated, in coordinates on $M$ and $N$
the tension field has components
\[
 \tau^\alpha(f)=\sum_{i,j=1}^m g^{ij}
 \left(\partial_i\partial_j f^\alpha
       -\sum_{a=1}^m\Gamma^a_{ij}\partial_a f^\alpha
       +\sum_{\beta,\gamma=1}^{\dim N}
         \Gamma^\alpha_{\beta\gamma}(f)
         \partial_i f^\beta\partial_j f^\gamma\right),
 \qquad \alpha\in\{1,\ldots,\dim N\}.
\]
On each fixed coordinate domain $\Omega\subseteq\mathbb R^m$,
the first term belongs to $H^{s-2}(\Omega)$; the coefficients on
the base are smooth. The first derivatives of the components of $f$
belong to $H^{s-1}(\Omega)$, which is an algebra, and
$\Gamma^\alpha_{\beta\gamma}(f)$ depends smoothly on the components
of $f$ with values in $H^s(\Omega)$. Estimate \eqref{eq:modulo-sobolev-tension} controls
multiplication at the lower level. Each differentiation with respect to $f$
replaces one factor by its variation or differentiates a superposition;
it does not add a spatial derivative to the variation vectors. This proves
smoothness with values in the fibers $H^{s-2}(M,f^*TN)$. For the energy, the same
products are applied to $|df|^2$, followed by the continuous integral on $L^1(M)$. Hence
\[
 \mathcal E\colon H^s(M,N)\longrightarrow\mathbb R
\]
is smooth. Similarly, the tension field defines the smooth section
\begin{equation}
\label{eq:perdida-mapas-armonicos}
 \boldsymbol{\tau}\colon
 H^s(M,N)\longrightarrow\boldsymbol{\mathcal T}^{s-2},
 \qquad
 f\longmapsto(f,\boldsymbol{\tau}(f)),
 \qquad
 \pi_{\boldsymbol{\mathcal T}^{s-2}}\circ\boldsymbol{\tau}
 =
 \operatorname{id}_{H^s(M,N)}.
\end{equation}
To compute its vertical derivative at $f$, we use the trivialization
by parallel transport centered precisely at $f$. This gives the
continuous linear operator
\[
 D\boldsymbol{\tau}(f)\colon
 H^s(M,f^*TN)
 \longrightarrow
 H^{s-2}(M,f^*TN),
\]
whose expression on smooth maps is
\eqref{eq:linealizacion-tension-mapa}. If \(f\) is harmonic,
\(J_f=-D\boldsymbol{\tau}(f)\) has the same domain and codomain.
In another trivialization, differentiating the change of frame produces a term
multiplied by $\boldsymbol\tau(f)$; this is why the linearization at a
zero is independent of that choice. For nonharmonic maps, the
intrinsic expression above is the covariant derivative, not the ordinary
derivative of all representations simultaneously.

Thus \(-\boldsymbol{\tau}(f)\) does not generally belong to the tangent space
\(T_fH^s(M,N)\), but to the fiber of lower regularity
\(H^{s-2}(M,f^*TN)\). An alternative is to use a strong metric
on the mapping space. The representative of the differential with respect to
such a metric is obtained by formally regularizing through
\[
 (1+\nabla^*\nabla)^{-s}(-\boldsymbol{\tau}(f)).
\]
Its construction requires studying how the elliptic operators involved
depend on \(f\).

Finiteness of the energy requires only that $df$ be square integrable, whereas the manifold
structure and pointwise evaluation require regularity greater than
\(\frac{m}{2}\). This difference between the energy level and the geometric
level also appears in compactness questions. In dimension two,
for example, a sequence with bounded energy may concentrate energy and
produce bubbles; on a noncompact base, it may also lose energy at
infinity.

\section{Gradient descent and geodesics}
\label{sec:descenso-gradiente-geodesicas}

The geodesic equation can be studied as an equation for critical points of
energy on a manifold of curves. Three objects must be distinguished: the
geodesic, with parameter $t$; the descent trajectory, with parameter $s$;
and the gradient of energy, which depends on the metric chosen on the curve
space. A descent trajectory is generally neither a geodesic of that space
nor a geodesic of the target manifold. Our aim is to obtain a geodesic as a
limit of curves with decreasing energy.

We use the Hilbert framework of $H^1(J,N)$, where $J$ is the interval $[0,1]$ or the circle $S^1$ and $N$ is the target manifold, which combines a well-defined
gradient with an appropriate compactness condition. This scheme appears in
Moore \cite[Sections 1.9--1.11 and 2.1--2.4]{Moore2017GlobalAnalysis}.
We develop it from the calculus, mapping manifolds, and flow theory already
established. We also include a discrete version with step selection and
specify the additional hypotheses that give convergence of the entire
trajectory.

\subsection{Descent flow on a strong Hilbert manifold}
\label{subsec:flujo-descenso-hilbert}

\begin{proposition}[Gradient and descent direction]

\label{prop:gradiente-fuerte-descenso}

Let $(\mathcal M,\mathbf G)$ be a Hilbert manifold with a smooth strong
Riemannian metric, and let $F\colon\mathcal M\to\mathbb R$ be of class
$C^r$, $r\geq2$. There is a unique $C^{r-1}$ vector field
$\operatorname{grad}_{\mathbf G}F$ such that

\begin{equation}\label{eq:def-gradiente-fuerte-descenso}
 \mathbf G_x(\operatorname{grad}_{\mathbf G}F(x),\mathbf V)=dF_x(\mathbf V),
 \qquad \mathbf V\in T_x\mathcal M.
\end{equation}
\[
 \|\operatorname{grad}_{\mathbf G}F(x)\|_{\mathbf G}=\|dF_x\|_{\mathbf G'}.
\]

If $dF_x\ne0$, the minimum of $dF_x(\mathbf V)$ over unit vectors is
attained uniquely at
$-\operatorname{grad}_{\mathbf G}F(x)/\|\operatorname{grad}_{\mathbf G}F(x)\|_{\mathbf G}$.
Consequently, the vector field $-\operatorname{grad}_{\mathbf G}F$ has a
unique local flow.

\end{proposition}

\begin{proof}

By Proposition~\ref{prop:caracterizaciones-metrica-fuerte}, the musical
map $\flat_{\mathbf G}\colon T\mathcal M\to T'\mathcal M$ is a smooth
bundle isomorphism. The required field is $\flat_{\mathbf G}^{-1}\circ dF$.
Riesz representation proves the norm identity. Cauchy--Schwarz gives
$dF_x(\mathbf V)\geq-\|\operatorname{grad}_{\mathbf G}F(x)\|_{\mathbf G}$
for $\|\mathbf V\|_{\mathbf G}=1$, and its equality case gives the stated
direction. The field is $C^1$, hence locally Lipschitz; apply
Theorem~\ref{teo:existencia-curvas-integrales-banach}.

\end{proof}

\begin{theorem}[Global existence and subsequential convergence]

\label{teo:descenso-global-palais-smale}

Suppose that the length distance $d_{\mathbf G}$ is complete on each
component of $\mathcal M$ and that $F\in C^2(\mathcal M)$ is bounded below.
For every $x_0\in\mathcal M$, the problem

\begin{equation}\label{eq:flujo-gradiente-negativo}
 \frac{dx}{ds}(s)=-\operatorname{grad}_{\mathbf G}F(x(s)),\qquad x(0)=x_0,
\end{equation}

has a solution for every $s\geq0$. For $0\leq a\leq b$, it satisfies

\begin{align}
 F(x(a))-F(x(b))
 &=\int_a^b\|\operatorname{grad}_{\mathbf G}F(x(s))\|_{\mathbf G}^2\,ds,
 \label{eq:identidad-disipacion-gradiente}\\
 d_{\mathbf G}(x(a),x(b))^2
 &\leq(b-a)\bigl(F(x(a))-F(x(b))\bigr).
 \label{eq:control-distancia-descenso}
\end{align}

If, in addition, every sequence $(x_n)$ with bounded $(F(x_n))$ and
$\|dF_{x_n}\|_{\mathbf G'}\to0$ has a convergent subsequence in
$\mathcal M$, there are $s_n\to\infty$ and a critical point $x_*$ such
that $x(s_n)\to x_*$. This is the Palais--Smale condition with respect to
$\mathbf G$. Global existence does not require this condition.

\end{theorem}

\begin{proof}

The chain rule and \eqref{eq:def-gradiente-fuerte-descenso} give
$\frac{d}{ds}F(x(s))=-\|\operatorname{grad}_{\mathbf G}F(x(s))\|_{\mathbf G}^2$.
Integration yields \eqref{eq:identidad-disipacion-gradiente}. Distance is
bounded by the length of the trajectory between $a$ and $b$;
Cauchy--Schwarz and the preceding identity give
\eqref{eq:control-distancia-descenso}.

Let $[0,T)$ be the maximal forward interval. If $T<\infty$,
monotonicity and the lower bound imply that $F(x(s))$ has a limit as
$s\to T^{-}$. The distance inequality shows that $x(s)$ is Cauchy as
$s\to T^{-}$, hence converges to a point $y$. The strong distance induces
the manifold topology. In a chart around $y$, the integral equation for
$x$ passes to the limit at $T$, since the vector field is continuous and
locally bounded. Setting $x(T)=y$ and adjoining the local solution from
$y$ extends the solution past $T$, contradicting maximality. Thus
$T=\infty$.

Finally, \eqref{eq:identidad-disipacion-gradiente} implies
$\displaystyle \int_0^\infty\|\operatorname{grad}_{\mathbf G}F(x(s))\|_{\mathbf G}^2ds<\infty$.
For every integer $n$, there is $s_n\geq n$ with gradient norm less than
$1/n$; otherwise the integral on a half-line would be infinite.
Palais--Smale supplies a convergent subsequence. Continuity of $dF$ makes
its limit critical. The conclusion concerns a temporal subsequence;
the next proposition gives conditions for convergence of the entire path.

\end{proof}

\begin{proposition}[A criterion for convergence of the whole trajectory]

\label{prop:convergencia-total-desigualdad-gradiente}

Under the preceding theorem's assumptions, suppose that
$x(s_n)\to x_*$ for some $s_n\to\infty$ and that, in a neighborhood $U$
of $x_*$, there are $C>0$ and $0<\theta\leq1/2$ such that

\begin{equation}\label{eq:desigualdad-gradiente-local}
 |F(y)-F(x_*)|^{1-\theta}
 \leq C\|\operatorname{grad}_{\mathbf G}F(y)\|_{\mathbf G},\qquad y\in U.
\end{equation}

Then $x(s)\to x_*$ as $s\to\infty$, and the trajectory has a finite-length
tail. For sufficiently large $s$,

\[
 d_{\mathbf G}(x(s),x_*)\leq
 \frac C\theta\bigl(F(x(s))-F(x_*)\bigr)^\theta.
\]

If $\theta=1/2$, both the energy difference and the distance decay
exponentially, possibly with different constants.

\end{proposition}

\begin{proof}

Put $e(s)=F(x(s))-F(x_*)$. The function $e$ is nonincreasing, and the given subsequence implies
$\displaystyle\lim_{s\to\infty}e(s)=0$. If $e$ vanishes at some time, the dissipation identity forces
the trajectory to be constant thereafter. The given subsequence identifies
that constant value with $x_*$, and all conclusions hold on that tail.
Thus assume $e(s)>0$ for every $s\geq0$. While $e(s)>0$ and
$x(s)\in U$,

\[
 -\frac{d}{ds}e(s)^\theta
 =\theta e(s)^{\theta-1}\|\operatorname{grad}_{\mathbf G}F(x(s))\|_{\mathbf G}^2
 \geq\frac\theta C\|x'(s)\|_{\mathbf G}.
\]

Choose $\rho>0$ with $B_{d_{\mathbf G}}(x_*,\rho)\subset U$ and an
index $n_0$ such that, setting $s_0=s_{n_0}$,
\[
 d_{\mathbf G}(x(s_0),x_*)<\rho/3,
 \qquad (C/\theta)e(s_0)^\theta<\rho/3.
\]
Define
\[
 \tau:=\inf\{s\geq s_0\mid d_{\mathbf G}(x(s),x_*)\geq\rho\},
 \qquad \inf\varnothing:=+\infty.
\]
If $\tau<\infty$, continuity of the distance implies $\tau>s_0$,
$d_{\mathbf G}(x(\tau),x_*)=\rho$, and $x(s)\in U$ for $s_0\leq s<\tau$.
For $s_0<b<\tau$, the preceding differential inequality gives
\[
 L_{\mathbf G}(x\restriction_{[s_0,b]})
 \leq\frac C\theta\bigl(e(s_0)^\theta-e(b)^\theta\bigr)
 \leq\frac C\theta e(s_0)^\theta<\rho/3.
\]
By the triangle inequality, letting $b\to\tau^{-}$ gives
\[
 \rho=d_{\mathbf G}(x(\tau),x_*)
 \leq d_{\mathbf G}(x(s_0),x_*)+\frac C\theta e(s_0)^\theta
 <2\rho/3,
\]
a contradiction. Hence $\tau=\infty$, and the entire tail stays in $U$.
For $s_0\leq a<b$, we therefore have
\[
 d_{\mathbf G}(x(a),x(b))
 \leq L_{\mathbf G}(x\restriction_{[a,b]})
 \leq\frac C\theta\bigl(e(a)^\theta-e(b)^\theta\bigr)
 \leq\frac C\theta e(a)^\theta.
\]
Since $e(a)\to0$, the trajectory is Cauchy as $a\to\infty$. Completeness
gives a limit, and the convergent subsequence identifies it with $x_*$.
Letting $b\to\infty$ in the length estimate shows that the tail has finite
length and gives
\[
 d_{\mathbf G}(x(a),x_*)\leq\frac C\theta e(a)^\theta,
 \qquad a\geq s_0.
\]
If $\theta=1/2$, then $e'(s)\leq-C^{-2}e(s)$ throughout this tail, and
$e(s)\leq e(s_0)e^{-(s-s_0)/C^2}$. The distance bound also proves
exponential convergence of $x(s)$.

\end{proof}

\begin{corollary}[Nondegenerate critical point]

\label{cor:descenso-limite-no-degenerado}

Inequality \eqref{eq:desigualdad-gradiente-local}, with $\theta=1/2$,
holds if the Hessian of $F$ at $x_*$ induces a bounded isomorphism
$T_{x_*}\mathcal M\to T'_{x_*}\mathcal M$. Thus a subsequence converging
to such a point determines exponential convergence of the whole trajectory.

\end{corollary}

\begin{proof}
In a chart centered at $x_*$, put $f=F\circ\varphi^{-1}$,
$A=D^2f(0):H\to H'$, and $c_A=\|A^{-1}\|^{-1}>0$.
Continuity of $D^2f$ gives $\rho>0$ such that $B_H(0,\rho)$ lies
in the coordinate domain and

\[
 \|D^2f(z)-A\|_{\mathcal L^2(H;\mathbb R)}\leq c_A/2,
 \qquad \|z\|_H<\rho.
\]
Since $Df(0)=0$, integral Taylor expansion gives, for $\|v\|_H<\rho$,

\begin{align*}
 Df(v)-Av&=\int_0^1(D^2f(tv)-A)[v,\,\cdot\,]\,dt,\\
 \|Df(v)-Av\|_{H'}&\leq(c_A/2)\|v\|_H,\\
 \|Df(v)\|_{H'}&\geq(c_A/2)\|v\|_H,\\
 |f(v)-f(0)|&\leq C_1\|v\|_H^2,
 \qquad C_1=\tfrac12(\|A\|+c_A/2).
\end{align*}
The strong metric is continuous in this chart. After reducing $\rho$,
there is $C_g\geq1$ such that

\[
 C_g^{-1}\|w\|_H\leq
 \|d\varphi^{-1}_v w\|_{\mathbf G_{\varphi^{-1}(v)}}
 \leq C_g\|w\|_H,\qquad \|v\|_H<\rho,\quad w\in H.
\]
Consequently, $\|Df(v)\|_{H'}\leq
C_g\|dF_{\varphi^{-1}(v)}\|_{\mathbf G'}$. By Riesz and the preceding bounds,

\[
 |F(x)-F(x_*)|^{1/2}
 \leq\frac{2C_g\sqrt{C_1}}{c_A}\|\operatorname{grad}_{\mathbf G}F(x)\|_{\mathbf G_x},
 \qquad x\in\varphi^{-1}(B_H(0,\rho)).
\]
In infinite dimensions, zero kernel does not imply the isomorphism
hypothesis: the operator $(v_j)\mapsto(v_j/j)$ on $\ell^2(\mathbb N)$ is injective
and self-adjoint, but its image is not all of $\ell^2(\mathbb N)$.

\end{proof}

\subsection{The curve manifold and its first-order metric}
\label{subsec:curvas-h1-metrica-fuerte}

From now on, $(N,\mathbf h)$ is a smooth connected finite-dimensional
Riemannian manifold without boundary. We use the parameter spaces
$I=[0,1]$ and $S^1=\mathbb R/\mathbb Z$, with measure $dt$.
For $p,q\in N$, consider

\[
 \mathcal P_{p,q}(N)=\{\gamma\in H^1(I,N)\mid\gamma(0)=p,\ \gamma(1)=q\},
 \qquad \Lambda(N)=H^1(S^1,N).
\]

An $H^1(J,N)$ map is identified with its continuous representative. This is
possible because $1>\dim(I)/2$. Write $J$ for the parameter space and
$\mathcal P$ for either of the preceding manifolds. Endpoint conditions
for a closed curve are not imposed through an arbitrary trivialization:
all identifications on $S^1$ are those of the bundle over the circle.

Example~\ref{ej:curvas-esfera-extremos-banach} constructed
$\mathcal P_{p,q}(S^{d-1})$ using explicit charts. For a general target
manifold, use the exponential charts in
Theorem~\ref{teo:cartas-exponenciales-mapeos-sobolev} and study both
the curve manifold structure and its first-order metric and energy.

\begin{lemma}[Elementary curve estimates]

\label{lem:h1-curvas-estimaciones}

For $u\in H^1(J,\mathbb R^a)$,

\[
 \|u\|_{C^0}\leq\|u\|_{L^2(J,\mathbb R^a)}+\|u'\|_{L^2(J,\mathbb R^a)}
 \leq\sqrt2\|u\|_{H^1(J,\mathbb R^a)},\qquad
 |u(t)-u(\tau)|\leq|t-\tau|^{1/2}\|u'\|_{L^2(J,\mathbb R^a)}.
\]

On the circle, the second inequality is read along any parametrized arc
of length at most one. If $\gamma\in H^1(J,N)$, an $H^1(J,\gamma^*TN)$ section of
$\gamma^*TN$ has a weak covariant derivative $D_t\mathbf V\in L^2(J,\gamma^*TN)$, and

\begin{equation}\label{eq:cota-sup-campo-curva}
 \|\mathbf V\|_{C^0}
 \leq \|\mathbf V\|_{L^2(J,\gamma^*TN)}+\|D_t\mathbf V\|_{L^2(J,\gamma^*TN)}
 \leq\sqrt2\left(\int_J(|\mathbf V|_{\mathbf h}^2+
                  |D_t\mathbf V|_{\mathbf h}^2)\,dt\right)^{1/2}.
\end{equation}

\end{lemma}

\begin{proof}
For $J=I$, Lemma~\ref{lem:algebra-h1-intervalo}, applied to each
component, gives the absolutely continuous representative of $u$ and
the identity
\[
 u(t)-u(\tau)=\int_\tau^t u'(v)\,dv.
\]
On the circle, apply the same lemma in coordinate arcs, using
affine reparametrizations of the intervals. The continuous
representatives agree on intersections because they represent the same
class almost everywhere. Their lift to $[0,1]$ is absolutely
continuous by joining finitely many segments and has equal endpoint
values. The preceding identity applies along the lift of any arc
of length at most one. In both cases,
Cauchy--Schwarz gives
\[
 |u(t)-u(\tau)|
 \leq\int_{[t,\tau]}|u'(v)|\,dv
 \leq |t-\tau|^{\frac{1}{2}}\|u'\|_{L^2(J,\mathbb R^a)},
\]
where $[t,\tau]$ denotes the segment between the chosen parameters.

To deduce the vector uniform bound from the scalar estimate in
Lemma~\ref{lem:algebra-h1-intervalo}, fix $\varepsilon>0$ and set
\[
 w_\varepsilon(t):=
 \bigl(|u(t)|^2+\varepsilon^2\bigr)^{\frac{1}{2}}.
\]
The product proved in that lemma and the inverse square root in
Example~\ref{ej:funciones-elementales-algebras} show that
$w_\varepsilon\in H^1(J,\mathbb R)$: the square root of a positive
function is its product with its inverse square root.
The same rules give, almost everywhere,
\[
 w_\varepsilon'
 =\frac{u\cdot u'}{w_\varepsilon},
 \qquad |w_\varepsilon'|\leq |u'|,
 \qquad
 \|w_\varepsilon\|_{L^2(J)}^2
 =\|u\|_{L^2(J,\mathbb R^a)}^2+\varepsilon^2.
\]
On the circle, apply these rules to the periodic lift;
the matching endpoint values cancel the boundary term in
integration by parts. Since $J$ has measure one,
the scalar estimate already proved implies
\[
 \|w_\varepsilon\|_{C^0}
 \leq\|w_\varepsilon\|_{L^2(J)}
       +\|w_\varepsilon'\|_{L^2(J)}
 \leq\bigl(\|u\|_{L^2(J,\mathbb R^a)}^2+\varepsilon^2\bigr)^{\frac{1}{2}}
       +\|u'\|_{L^2(J,\mathbb R^a)}.
\]
The inequality $0\leq w_\varepsilon-|u|\leq\varepsilon$ allows a
uniform passage to the limit. This gives the first bound in the statement;
the inequality $A+B\leq\sqrt2(A^2+B^2)^{\displaystyle\frac{1}{2}}$ for
$A,B\geq0$ gives its second expression.

Now consider a section $\mathbf V\in H^1(J,\gamma^*TN)$.
Continuity of $\gamma$ and compactness of $J$ give finitely many
segments whose images lie in compact subsets of charts of $N$.
On each segment, the components of $\gamma$ and $\mathbf V$ are $H^1$
and have the preceding continuous representatives. Write
\[
 (D_t\mathbf V)^a=(V^a)'
 +\sum_{b,c=1}^{\dim N}
       \Gamma^a_{bc}(\gamma)(\gamma^b)'V^c.
\]
The coefficients $\Gamma^a_{bc}(\gamma)$ and $V^c$ are bounded on
the segment, whereas $(\gamma^b)'$ and $(V^a)'$ belong to $L^2$.
Thus each component of this expression belongs to $L^2$.
The weak product rule in Lemma~\ref{lem:algebra-h1-intervalo}
and the coordinate chain rule in
Lemma~\ref{lem:superposicion-sobolev-geometrica} give the
transformation law for these components. The local expressions
therefore define a global section $D_t\mathbf V\in L^2(J,\gamma^*TN)$.

To verify the intrinsic estimate, define
\[
 z_\varepsilon(t):=
 \bigl(|\mathbf V(t)|_{\mathbf h}^2+\varepsilon^2\bigr)^{\frac{1}{2}}.
\]
The components of $\mathbf h(\gamma)$ are $H^1$ by the same
superposition lemma. The preceding products and inverse square root show
that $z_\varepsilon$ is $H^1$ on each segment. These scalar functions
agree on intersections and join to a function in $H^1(J)$.
Metric compatibility, written as
\[
 \partial_c h_{ab}
 =\sum_{d=1}^{\dim N}
       \bigl(\Gamma^d_{ca}h_{db}+\Gamma^d_{cb}h_{ad}\bigr),
\]
and the weak product and chain rules give
\[
 \bigl(|\mathbf V|_{\mathbf h}^2\bigr)'
 =2\mathbf h(D_t\mathbf V,\mathbf V),
 \qquad
 z_\varepsilon'
 =\frac{\mathbf h(D_t\mathbf V,\mathbf V)}{z_\varepsilon}
\]
almost everywhere. Cauchy--Schwarz in each fiber gives
\[
 |z_\varepsilon'|
 \leq |D_t\mathbf V|_{\mathbf h}
       \frac{|\mathbf V|_{\mathbf h}}
            {\bigl(|\mathbf V|_{\mathbf h}^2+\varepsilon^2\bigr)^{\frac{1}{2}}}
 \leq |D_t\mathbf V|_{\mathbf h}.
\]
Consequently, the same scalar estimate implies
\[
 \|z_\varepsilon\|_{C^0}
 \leq
 \bigl(\|\mathbf V\|_{L^2(J,\gamma^*TN)}^2+\varepsilon^2\bigr)^{\frac{1}{2}}
 +\|D_t\mathbf V\|_{L^2(J,\gamma^*TN)}.
\]
Since $0\leq z_\varepsilon-|\mathbf V|_{\mathbf h}\leq\varepsilon$,
letting $\varepsilon\to0^+$ proves the first inequality in
\eqref{eq:cota-sup-campo-curva}. The inequality for $A+B$ used above
proves the second.\qedhere
\end{proof}

\begin{proposition}[Hilbert structure and energy]

\label{prop:curvas-h1-energia-metrica}

The spaces $\mathcal P_{p,q}(N)$ and $\Lambda(N)$ are smooth Hilbert
manifolds. Their tangent spaces are, respectively,

\begin{align*}
 T_\gamma\mathcal P_{p,q}(N)
 &=\{\mathbf V\in H^1(J,\gamma^*TN)\mid\mathbf V(0)=0,\ \mathbf V(1)=0\},\\
 T_\gamma\Lambda(N)&=H^1(S^1,\gamma^*TN).
\end{align*}
\begin{equation}\label{eq:metrica-h1-curvas-intrinseca}
 \mathbf G_\gamma(\mathbf V,\mathbf W)
 =\int_J\bigl(\mathbf h(\mathbf V,\mathbf W)+
              \mathbf h(D_t\mathbf V,D_t\mathbf W)\bigr)\,dt
\end{equation}

is a smooth strong Riemannian metric. The energy

\begin{equation}\label{eq:energia-curvas-h1}
 \mathcal E(\gamma)=\frac12\int_J|\dot\gamma(t)|_{\mathbf h}^2\,dt
\end{equation}

is smooth and satisfies

\begin{align}
 d\mathcal E_\gamma(\mathbf V)
 &=\int_J\mathbf h(\dot\gamma,D_t\mathbf V)\,dt,
 \label{eq:primera-variacion-curvas-h1}\\
 \|d\mathcal E_\gamma\|_{\mathbf G'}
 &\leq\sqrt{2\mathcal E(\gamma)}.
 \label{eq:cota-gradiente-energia-curvas}
\end{align}

\end{proposition}

\begin{proof}

Apply Theorem~\ref{teo:cartas-exponenciales-mapeos-sobolev}
with $r=1$, $p=2$, and a one-dimensional source. The chart centered at
$\gamma\in H^1(J,N)$ is
$\operatorname{Exp}_\gamma(\mathbf V)(t)
=\exp_{\gamma(t)}\mathbf V(t)$ on a neighborhood of the zero section.
On the interval, if $\gamma(0)=p$ and $\gamma(1)=q$, injectivity
of the exponential maps on the chart domains shows that
$\operatorname{Exp}_\gamma(\mathbf V)$ has the same endpoints
if and only if $\mathbf V(0)=\mathbf V(1)=0$.

The evaluation operator
\[
 L_\gamma:H^1(I,\gamma^*TN)\longrightarrow T_pN\oplus T_qN,
 \qquad L_\gamma(\mathbf V)=(\mathbf V(0),\mathbf V(1)),
\]
is linear and continuous by Lemma~\ref{lem:h1-curvas-estimaciones}.
Set $n:=\dim N$. Construct a right inverse explicitly. Choose smooth frames
near $p$ and $q$, orthonormal at those points, with endpoint bases
$(e_j)_{j=1}^{n}$ and $(f_j)_{j=1}^{n}$. Continuity
of $\gamma$ allows their restriction to disjoint intervals near
$0$ and $1$. Denote the induced fields along
$\gamma$ by $e_j^0(t)$ and $f_j^1(t)$, so that
$e_j^0(0)=e_j$ and $f_j^1(1)=f_j$. Choose smooth functions
$\chi_0,\chi_1$ supported in those intervals and equal to one
near the corresponding endpoint, and set
\[
 \boldsymbol\eta_j(t)=\chi_0(t)e_j^0(t),
 \qquad
 \boldsymbol\zeta_j(t)=\chi_1(t)f_j^1(t).
\]
Extend them by zero outside their intervals. They are $H^1$ sections:
in a trivialization their components are products of smooth
functions of $t$ with smooth superpositions of $\gamma$; the cutoff
supports prevent jumps at the interior endpoints of the segments.
For $v=\displaystyle\sum_{j=1}^{n}v^je_j$ and $w=\displaystyle\sum_{j=1}^{n}w^jf_j$, define
\[
 R_\gamma(v,w)
 :=\sum_{j=1}^{n}v^j\boldsymbol\eta_j
   +\sum_{j=1}^{n}w^j\boldsymbol\zeta_j.
\]
Then $L_\gamma R_\gamma=\operatorname{id}$. Moreover,
\[
 \|R_\gamma(v,w)\|_{H^1}
 \leq C_\gamma
       \bigl(|v|_{\mathbf h}^2+|w|_{\mathbf h}^2\bigr)^{\frac{1}{2}},
 \qquad
 C_\gamma:=
 \left(\sum_{j=1}^{n}
    \bigl(\|\boldsymbol\eta_j\|_{H^1}^2
         +\|\boldsymbol\zeta_j\|_{H^1}^2\bigr)\right)^{\frac{1}{2}},
\]
by the triangle inequality and Cauchy--Schwarz. Thus,
\[
 H^1(I,\gamma^*TN)
 =\ker L_\gamma\oplus R_\gamma(T_pN\oplus T_qN),
 \qquad
 P_\gamma=\operatorname{id}-R_\gamma L_\gamma
\]
is the decomposition given by a continuous projection onto the kernel.
Restricting the exponential chart to that kernel gives the Hilbert
structure on $\mathcal P_{p,q}(N)$ and the stated tangent
identification. For $\Lambda(N)$, apply the same theorem
directly with source $S^1$, without choosing additional endpoint identifications.

Expression \eqref{eq:metrica-h1-curvas-intrinseca} is positive definite.
In finitely many trivializations,
$\|D_t\mathbf V\|_{L^2(J,\gamma^*TN)}\leq C(\|\mathbf V'\|_{L^2(J,\mathbb R^n)}+
\|\dot\gamma\|_{L^2(J,\gamma^*TN)}\|\mathbf V\|_{C^0})$, giving one comparison with
the $H^1(J,\mathbb R^n)$ norm. Solving for $\mathbf V'$ and using
\eqref{eq:cota-sup-campo-curva} gives the reverse comparison. The
constants are local in $\gamma$ and uniform on an $H^1(J,\gamma^*TN)$ neighborhood.
Thus the metric is strong.

For smoothness, in coordinates each integrand is a finite sum of products
of smooth coefficients evaluated at the local representations of
$\gamma\in H^1(J,N)$, components of fields in $H^1(J,\gamma^*TN)$,
and first derivatives of those representations and components.
These derivatives belong to $L^2$ on each coordinate subinterval.
Smooth superposition from Lemma~\ref{lem:superposicion-sobolev-geometrica},
the continuous products $H^1(J)\times H^1(J)\to H^1(J)$ and
$H^1(J)\times L^2(J)\to L^2(J)$ provided by
Lemma~\ref{lem:algebra-h1-intervalo} on each coordinate segment, and the
integrated inner product $L^2(J)\times L^2(J)\to\mathbb R$ prove smoothness in
all variables. The same argument applies to $\mathcal E$.
Differentiating $\frac12\displaystyle\sum_{a,b=1}^{\dim N}h_{ab}(x)(x^a)'(x^b)'$ gives

\[
 \int_J\left(
 \displaystyle\sum_{a,b=1}^{\dim N}h_{ab}(x)(x^a)'(V^b)'
 +\frac12\displaystyle\sum_{a,b,c=1}^{\dim N}
       \partial_c h_{ab}(x)V^c(x^a)'(x^b)'
 \right)dt.
\]

The Christoffel-symbol formula turns this into
\eqref{eq:primera-variacion-curvas-h1}. When several charts are needed,
subdivide $J$ and add the integrals; no integration by parts has been
used, so no artificial terms arise at subdivision points.
Finally, Cauchy--Schwarz gives
$|d\mathcal E_\gamma(\mathbf V)|\leq
\|\dot\gamma\|_{L^2(J,\gamma^*TN)}\|D_t\mathbf V\|_{L^2(J,\gamma^*TN)}$, proving the last bound.\qedhere

\end{proof}

\begin{proposition}[Comparison with the metric induced by an immersion]

\label{prop:comparacion-metricas-h1-inmersion}

Let $\iota\colon(N,\mathbf h)\to\mathbb R^L$ be a smooth isometric
immersion, with second fundamental form $\mathrm{II}$.
The metric induced by the usual inner product on
$H^1(J,\mathbb R^L)$ is

\begin{equation}\label{eq:metrica-h1-curvas-extrinseca}
 \mathbf G^\iota_\gamma(\mathbf V,\mathbf W)
 =\int_J\left(
 \langle d\iota_\gamma\mathbf V,d\iota_\gamma\mathbf W\rangle
 +\left\langle\frac{d}{dt}(d\iota_\gamma\mathbf V),
                  \frac{d}{dt}(d\iota_\gamma\mathbf W)\right\rangle
 \right)dt.
\end{equation}

Here derivatives are weak and the inner product is Euclidean.
For an embedding with $N$ compact and $J=S^1$, this is the construction
in \cite[Example 1.9.1]{Moore2017GlobalAnalysis}.
For every $\gamma\in\mathcal P$ and all
$\mathbf V,\mathbf W\in T_\gamma\mathcal P$, we have

\begin{equation}\label{eq:comparacion-gauss-metricas-h1}
 \mathbf G^\iota_\gamma(\mathbf V,\mathbf W)
 =\mathbf G_\gamma(\mathbf V,\mathbf W)
  +\int_J\left\langle\mathrm{II}_\gamma(\dot\gamma,\mathbf V),
                        \mathrm{II}_\gamma(\dot\gamma,\mathbf W)
             \right\rangle dt.
\end{equation}

Let $K\subset N$ be compact and let $0\leq a<\infty$. Write

\[
 B_K=\sup_{y\in K}\|\mathrm{II}_y\|,\qquad
 C_{K,a}=\sqrt{1+4aB_K^2},
\]

where the norm of $\mathrm{II}_y$ is the bilinear operator norm for
$\mathbf h_y$ and the Euclidean inner product. If
$\gamma(J)\subset K$ and $\mathcal E(\gamma)\leq a$, then

\begin{align}
 \|\mathbf V\|_{\mathbf G_\gamma}
 &\leq\|\mathbf V\|_{\mathbf G^\iota_\gamma}
 \leq C_{K,a}\|\mathbf V\|_{\mathbf G_\gamma},
 \label{eq:equivalencia-metricas-h1-inmersion}\\
 C_{K,a}^{-1}\|\xi\|_{\mathbf G_\gamma'}
 &\leq\|\xi\|_{(\mathbf G^\iota_\gamma)'}
 \leq\|\xi\|_{\mathbf G_\gamma'},
 \qquad \xi\in(T_\gamma\mathcal P)'.
 \label{eq:equivalencia-duales-h1-inmersion}
\end{align}

In particular, along a sequence satisfying these two uniform bounds,
the condition that the norm of the differential tend to zero is the
same for both metrics.

\end{proposition}

\begin{proof}

The weak chain rule in coordinates and the Gauss formula give, for
almost every $t\in J$,

\[
 \frac{d}{dt}(d\iota_\gamma\mathbf V)
 =d\iota_\gamma(D_t\mathbf V)
   +\mathrm{II}_\gamma(\dot\gamma,\mathbf V).
\]

Both terms on the right belong to $L^2(J,\mathbb R^L)$:
the first because $\mathbf V$ is $H^1(J,\gamma^*TN)$, and the second because
$\dot\gamma\in L^2(J,\gamma^*TN)$, $\mathbf V$ is continuous, and the coefficients
of $\mathrm{II}$ are bounded on the compact image of $\gamma$.
The first term is tangent to the immersion and the second is normal.
Their orthogonality and the isometric property of $\iota$ prove
\eqref{eq:comparacion-gauss-metricas-h1}. Moreover,

\[
 \int_J|\mathrm{II}_\gamma(\dot\gamma,\mathbf V)|^2dt
 \leq B_K^2\|\dot\gamma\|_{L^2(J,\gamma^*TN)}^2
                    \|\mathbf V\|_{C^0(J)}^2
 \leq4aB_K^2\|\mathbf V\|_{\mathbf G_\gamma}^2,
\]

by \eqref{eq:cota-sup-campo-curva}. This proves
\eqref{eq:equivalencia-metricas-h1-inmersion}. Taking suprema over the
corresponding unit balls gives
\eqref{eq:equivalencia-duales-h1-inmersion}. Applying these bounds to the
differential term by term, with the same $C_{K,a}$, proves the last
assertion.\qedhere

\end{proof}

\begin{lemma}[Exponential coordinates along a curve]

\label{lem:exp-curvas-estimaciones-uniformes}

Let $K\subset N$ be compact. There is $\delta>0$ such that, for
$y\in K$ and $|v|_{\mathbf h}<\delta$, the exponential map, its normal
inverse, and parallel transport along the radial segment are smoothly
defined. If $\gamma\in H^1(J,N)$ takes values in $K$,
$\mathbf u\in H^1(J,\gamma^*TN)$, and $\|\mathbf u\|_{C^0}<\delta$, write
$R_\gamma(\mathbf u)(t)=\exp_{\gamma(t)}\mathbf u(t)$.
There are smooth endomorphisms $A(y,v)$ and $B(y,v)$ of $T_yN$ such that

\begin{equation}\label{eq:velocidad-exp-curvas}
 P_{R_\gamma(\mathbf u)(t)\to\gamma(t)}
       \frac{d}{dt}R_\gamma(\mathbf u)(t)
 = A(\gamma(t),\mathbf u(t))D_t\mathbf u(t)
   +B(\gamma(t),\mathbf u(t))\dot\gamma(t).
\end{equation}
\[
 A(y,0)=B(y,0)=\operatorname{id}_{T_yN}.
\]

After reducing $\delta$, these coefficients and their derivatives of
any fixed order are uniformly bounded, and $A(y,v)$ has a uniformly
bounded inverse. Moreover, for every $0\leq a<\infty$ there are
$\rho_a>0$, $L_a<\infty$, and $C_a<\infty$ such that, if
$\mathcal E(\gamma)\leq a$, $\gamma(J)\subset K$, and
$\|\mathbf u\|_{\mathbf G_\gamma}\leq\rho_a$, then

\begin{align}
 \|D^2(\mathcal E\circ R_\gamma)(\mathbf u)\|
    _{\mathcal L^2(T_\gamma\mathcal P;\mathbb R)}&\leq L_a,
 \label{eq:taylor-uniforme-energia-retraccion}\\
 \|dR_\gamma(\mathbf u)\mathbf V\|_{\mathbf G_{R_\gamma(\mathbf u)}}
   &\leq C_a\|\mathbf V\|_{\mathbf G_\gamma}.
 \label{eq:diferencial-retraccion-curvas}
\end{align}

The operator norm in the first line uses $\mathbf G_\gamma$ in both
variables.

\end{lemma}

\begin{proof}

If $K=\varnothing$, the statement is vacuous. Suppose $K\neq\varnothing$.
For each $x\in K$, choose an open set $O_x\Subset N$ containing $x$.
Lemma~\ref{lem:vecindades-normales-convexas-locales} gives an open set
$W_x\subseteq O_x$ containing $x$ and $\varepsilon_x>0$ such that, for
every $y\in W_x$, $\exp_y$ is a diffeomorphism on
$B_{\mathbf h_y}(0,2\varepsilon_x)$, with image contained in $O_x$.
Choose $x_1,\ldots,x_m\in K$ such that
$K\subseteq\bigcup_{j=1}^mW_{x_j}$ and set
\[
 U:=\bigcup_{j=1}^mW_{x_j},\qquad
 \delta_0:=\min_{1\leq j\leq m}\varepsilon_{x_j}>0,\qquad
 \mathcal D:=\{(y,v)\in TN\mid y\in U,\ |v|_{\mathbf h}<2\delta_0\}.
\]
The map $F\colon\mathcal D\to N\times N$, defined by
$F(y,v):=(y,\exp_yv)$, is injective: its first component determines
$y$, and the exponential map is injective on each fiber ball. Its
differential is invertible, since the first component has identity
differential and the vertical differential of the second is invertible.
The inverse function theorem therefore makes $F$ a diffeomorphism onto
the open set $F(\mathcal D)$. Thus the normal inverse depends smoothly
on both endpoints on a common domain. Parallel transport along
$s\mapsto\exp_y(sv)$, $0\leq s\leq1$, is smooth in $(y,v)$ by smooth
dependence for the linear transport equation.

Use the connection to split $T_vTN\cong T_yN\oplus T_yN$ into horizontal
and vertical components. The velocity of $(\gamma(t),\mathbf u(t))$
has components $(\dot\gamma(t),D_t\mathbf u(t))$. Compose the two
components of $d\exp_v$ with inverse radial transport; these are
$B(y,v)$ and $A(y,v)$, respectively. At $v=0$, both are the identity.
The weak chain rule gives \eqref{eq:velocidad-exp-curvas} for
$\gamma\in H^1(J,N)$ and $\mathbf u\in H^1(J,\gamma^*TN)$,
for almost every $t\in J$.

The closed disk bundle
$\{(y,v)\mid y\in K,\ |v|_{\mathbf h}\leq\delta_0\}$ is compact and
contained in $\mathcal D$. In a finite atlas over it, smoothness bounds
the coefficients and their derivatives of every fixed order.
There exists $0<\delta<\delta_0$ such that
$\|A(y,v)-\operatorname{id}\|<1/2$ for $y\in K$ and $|v|_{\mathbf h}<\delta$.
Indeed, if no such radius existed, there would be $y_j\in K$ and
$v_j\in T_{y_j}N$ with $|v_j|_{\mathbf h}\to0$ and
$\|A(y_j,v_j)-\operatorname{id}\|\geq1/2$. A subsequence satisfies
$y_j\to y\in K$ and $(y_j,v_j)\to(y,0)$ in $TN$, contradicting
continuity of $A$ and $A(y,0)=\operatorname{id}$. The Neumann series
then gives $\|A(y,v)^{-1}\|\leq2$ on this domain.

We detail the estimates needed for the discrete method. Suppress
$(\gamma(t),\mathbf u(t))$ from the coefficients and put
$w=A D_t\mathbf u+B\dot\gamma$. Since parallel transport is isometric,

\[
 (\mathcal E\circ R_\gamma)(\mathbf u)=\frac12\int_J|w|_{\mathbf h}^2dt.
\]
\begin{align*}
 D w[\mathbf V]&= (D_vA[\mathbf V])D_t\mathbf u
                   +A D_t\mathbf V+(D_vB[\mathbf V])\dot\gamma,\\
 D^2w[\mathbf V,\mathbf W]
 &= (D_v^2A[\mathbf V,\mathbf W])D_t\mathbf u
    +(D_vA[\mathbf V])D_t\mathbf W
    +(D_vA[\mathbf W])D_t\mathbf V
    +(D_v^2B[\mathbf V,\mathbf W])\dot\gamma.
\end{align*}

Choose $\rho_a$ with $\sqrt2\rho_a<\delta/2$. Bound
\eqref{eq:cota-sup-campo-curva} places the whole ball of radius $\rho_a$
in the normal domain. The preceding formulas, coefficient bounds, and
$\|\dot\gamma\|_{L^2(J,\gamma^*TN)}\leq\sqrt{2a}$ give

\begin{align*}
 \|w\|_{L^2(J,\gamma^*TN)}&\leq C(\rho_a+\sqrt{2a}),\\
 \|Dw[\mathbf V]\|_{L^2(J,\gamma^*TN)}&\leq C(1+\rho_a+\sqrt{2a})\|\mathbf V\|_{\mathbf G_\gamma},\\
 \|D^2w[\mathbf V,\mathbf W]\|_{L^2(J,\gamma^*TN)}
 &\leq C(1+\rho_a+\sqrt{2a})
       \|\mathbf V\|_{\mathbf G_\gamma}\|\mathbf W\|_{\mathbf G_\gamma}.
\end{align*}

Differentiating the integral of $|w|^2/2$ twice gives
$\displaystyle \int_J\mathbf h(Dw[\mathbf V],Dw[\mathbf W])dt$ and
$\displaystyle \int_J\mathbf h(w,D^2w[\mathbf V,\mathbf W])dt$.
Cauchy--Schwarz proves \eqref{eq:taylor-uniforme-energia-retraccion}.
For \eqref{eq:diferencial-retraccion-curvas}, the field
$dR_\gamma(\mathbf u)\mathbf V$ is the image of $\mathbf V$ under the
vertical component of $d\exp$. Its covariant $t$-derivative is a sum of
one term bounded by $C|D_t\mathbf V|$ and terms bounded by
$C(|\dot\gamma|+|D_t\mathbf u|)|\mathbf V|$. Its $L^2(J,\gamma^*TN)$ norm, together
with that of the field, is at most $C_a\|\mathbf V\|_{\mathbf G_\gamma}$
by \eqref{eq:cota-sup-campo-curva}.\qedhere

\end{proof}

\begin{theorem}[Completeness of the curve-space metric]

\label{teo:completitud-metrica-h1-curvas}

If $(N,d_{\mathbf h})$ is complete, the length distance of
\eqref{eq:metrica-h1-curvas-intrinseca} is complete on every component
of $\mathcal P_{p,q}(N)$ and $\Lambda(N)$.

\end{theorem}

\begin{proof}

For a piecewise $C^1$ path $s\mapsto\alpha_s\in\mathcal P$, evaluation
at $t$ and \eqref{eq:cota-sup-campo-curva} give
$d_{\mathbf h}(\alpha_0(t),\alpha_1(t))\leq\sqrt2 L_{\mathbf G}(\alpha)$.
Taking the supremum in $t$ and the infimum over paths yields

\begin{equation}\label{eq:distancia-uniforme-controlada-h1}
 \sup_{t\in J}d_{\mathbf h}(\gamma(t),\eta(t))
       \leq\sqrt2 d_{\mathbf G}(\gamma,\eta).
\end{equation}

Moreover, $Q(\gamma)=\sqrt{1+2\mathcal E(\gamma)}$ satisfies
$|dQ_\gamma(\mathbf V)|\leq\|\mathbf V\|_{\mathbf G}$ by
\eqref{eq:cota-gradiente-energia-curvas}. It is therefore $1$--Lipschitz
for $d_{\mathbf G}$. Let $(\gamma_n)$ be Cauchy in one component.
The preceding inequality and completeness of $N$ give a continuous
uniform limit $\gamma\colon J\to N$ preserving the fixed endpoints
or periodicity. The energies of $\gamma_n$ are bounded, since
$Q(\gamma_n)$ is Cauchy.

Cover the compact image of $\gamma$ by finitely many relatively compact
charts and subdivide $J$ into segments mapped inside them. The same holds
for $\gamma_n$ when $n$ is large. On each segment $J_0\subseteq J$, their coordinate
expressions are bounded in $H^1(J_0,\mathbb R^n)$: compactness bounds the $L^2(J_0,\mathbb R^n)$ norm,
and energy together with equivalence to the Euclidean metric bounds
the derivative. Weak compactness of $H^1(J_0,\mathbb R^n)$ balls gives a weak limit;
uniform convergence identifies it with the coordinates of $\gamma$.
Thus $\gamma\in H^1(J,N)$.

Strong convergence remains to be proved; uniform convergence alone
is insufficient. In a normal tube around $\gamma$, write
$\Phi_\gamma(\eta)(t)=\exp_{\gamma(t)}^{-1}\eta(t)$.
If $\eta$ lies in that tube and $\mathcal E(\eta)\leq b$, then

\begin{equation}\label{eq:control-log-curvas}
 \|d\Phi_\gamma(\eta)\mathbf W\|_{H^1(J,\gamma^*TN)}
 \leq C_b\|\mathbf W\|_{\mathbf G_\eta}.
\end{equation}

Indeed, if $L(y,z)=\exp_y^{-1}z$, the field on the left is
$D_zL(\gamma,\eta)\mathbf W$. Its $t$-derivative contains
$D_zL\,D_t\mathbf W$ and terms involving
$\dot\gamma\,\mathbf W$ and $\dot\eta\,\mathbf W$.
All coefficients are bounded on a smaller closed tube over the compact
image of $\gamma$. Using \eqref{eq:cota-sup-campo-curva},
$\|\dot\eta\|_{L^2(J,\eta^*TN)}\leq\sqrt{2b}$, and $\gamma\in H^1(J,N)$ gives
\eqref{eq:control-log-curvas}.

Fix $\delta>0$ so that the preceding estimate holds whenever
$\displaystyle\sup_{t\in J}d_{\mathbf h}(\eta(t),\gamma(t))<\delta$.
Choose $n_0$ such that, for $n,m\geq n_0$,

\[
 \sup_{t\in J}d_{\mathbf h}(\gamma_n(t),\gamma(t))<\frac\delta4,
 \qquad
 d_{\mathbf G}(\gamma_n,\gamma_m)<\frac{\delta}{8\sqrt2}.
\]

For $0<\varepsilon<\delta/(8\sqrt2)$, the definition of length distance
provides a piecewise $C^1$ path $\alpha\colon[0,1]\to\mathcal P$,
with endpoints $\gamma_n,\gamma_m$, such that

\[
 L_{\mathbf G}(\alpha)
 <d_{\mathbf G}(\gamma_n,\gamma_m)+\varepsilon
 <\frac{\delta}{4\sqrt2}.
\]

Inequality \eqref{eq:distancia-uniforme-controlada-h1}, applied also to
each initial segment of $\alpha$, shows that, for
$0\leq s\leq1$,

\[
 \sup_{t\in J}d_{\mathbf h}(\alpha_s(t),\gamma(t))
 \leq\sqrt2 L_{\mathbf G}(\alpha)
       +\sup_{t\in J}d_{\mathbf h}(\gamma_n(t),\gamma(t))
 <\frac\delta2.
\]

Next, put
$\displaystyle Q_*:=\sup_{n\geq n_0}Q(\gamma_n)<\infty$.
The $1$--Lipschitz property of $Q$ gives

\[
 \mathcal E(\alpha_s)\leq b,
 \qquad
 b:=\frac12\left[\left(Q_*+\frac{\delta}{4\sqrt2}\right)^2-1\right].
\]

The same constant $C_b$ therefore applies along all these paths.
Integrate \eqref{eq:control-log-curvas} over $s\in[0,1]$ and take
the limit as $\varepsilon\to0^+$ to obtain

\[
 \|\Phi_\gamma(\gamma_n)-\Phi_\gamma(\gamma_m)\|_{H^1(J,\gamma^*TN)}
 \leq C_b\,d_{\mathbf G}(\gamma_n,\gamma_m).
\]

Thus $\mathbf u_n:=\Phi_\gamma(\gamma_n)$ is Cauchy in
$H^1(J,\gamma^*TN)$. Its uniform limit is zero; the continuous inclusion
into $C^0$ also identifies its $H^1(J,\gamma^*TN)$ limit with zero.
Continuity of the strong metric in the chart $R_\gamma$ allows us to
choose $r>0$ and $C_R<\infty$ such that

\[
 \|dR_\gamma(\mathbf v)\mathbf V\|_{\mathbf G_{R_\gamma(\mathbf v)}}
 \leq C_R\|\mathbf V\|_{H^1(J,\gamma^*TN)}
 \quad\text{if}\quad \|\mathbf v\|_{H^1(J,\gamma^*TN)}<r.
\]

For large $n$, the segment $s\mapsto R_\gamma(s\mathbf u_n)$,
$0\leq s\leq1$, lies in this chart and has length at most
$C_R\|\mathbf u_n\|_{H^1(J,\gamma^*TN)}$. It preserves the endpoint
or periodicity conditions. Thus $\gamma$ is in that component and
$d_{\mathbf G}(\gamma_n,\gamma)\to0$.\qedhere

\end{proof}

\subsection{The Sobolev gradient and the geodesic equation}
\label{subsec:gradiente-sobolev-geodesicas}

\begin{proposition}[Weak computation of the gradient]

\label{prop:gradiente-sobolev-energia}

Let $\mathcal H_\gamma=T_\gamma\mathcal P$ with inner product
$\mathbf G_\gamma$, and let $\mathcal H_\gamma'$ be its continuous dual.
The operator

\[
 P_\gamma\colon\mathcal H_\gamma\longrightarrow\mathcal H_\gamma',
 \qquad (P_\gamma\mathbf Z)(\mathbf V)=\mathbf G_\gamma(\mathbf Z,\mathbf V),
\]

is an isometric isomorphism. The gradient
$\mathbf Z_\gamma=\operatorname{grad}_{\mathbf G}\mathcal E(\gamma)$
is the unique solution of

\begin{equation}\label{eq:gradiente-geodesicas-forma-debil}
 \int_J\bigl(\mathbf h(\mathbf Z_\gamma,\mathbf V)+
           \mathbf h(D_t\mathbf Z_\gamma,D_t\mathbf V)\bigr)dt
 =\int_J\mathbf h(\dot\gamma,D_t\mathbf V)dt,
 \qquad\mathbf V\in\mathcal H_\gamma.
\end{equation}

If $\gamma$ is smooth, its gradient is smooth and satisfies

\begin{equation}\label{eq:gradiente-geodesicas-green}
 (\operatorname{id}-D_t^2)\mathbf Z_\gamma=-D_t\dot\gamma,
\end{equation}

with Dirichlet conditions $\mathbf Z_\gamma(0)=\mathbf Z_\gamma(1)=0$
for fixed endpoints, or as an equation of sections over $S^1$ in the
periodic case.

\end{proposition}

\begin{proof}

Riesz gives invertibility of $P_\gamma$, and
\eqref{eq:primera-variacion-curvas-h1} gives
$\mathbf Z_\gamma=P_\gamma^{-1}d\mathcal E_\gamma$.
For smooth $\gamma$, integration by parts against smooth interior test
fields gives \eqref{eq:gradiente-geodesicas-green} in distributions.
In a smooth frame, the principal term is $-\mathbf Z''$ and the other
terms are combinations of $\mathbf Z'$ and $\mathbf Z$ with smooth
coefficients. Since $\mathbf Z\in H^1(J,\mathbb R^n)$, solving for $\mathbf Z''$
gives an $L^2(J,\mathbb R^n)$ right-hand side. Thus $\mathbf Z\in H^2(J,\mathbb R^n)$.
Differentiate the equation and repeat to obtain $\mathbf Z\in H^k(J,\mathbb R^n)$
for every $k$, hence smoothness. The argument works up to the endpoints
in a frame along $I$; Dirichlet traces are already imposed by
$\mathcal H_\gamma$. On the circle, use charts crossing the point where
the interval is identified.

\end{proof}

\begin{remark}[Two different gradients]

\label{obs:gradiente-h1-no-calor}

At a merely $H^1(J,N)$ curve, $-D_t\dot\gamma$ means the functional on the
right of \eqref{eq:gradiente-geodesicas-forma-debil}, not an $H^1(J,\gamma^*TN)$
field. The inverse $P_\gamma^{-1}$ converts this covector into a tangent
vector. The informal equality ``gradient equals minus acceleration''
uses the weak $L^2(J,\gamma^*TN)$ pairing, not the strong metric
\eqref{eq:metrica-h1-curvas-intrinseca}. Our descent is therefore

\[
 \partial_s\gamma=-P_\gamma^{-1}d\mathcal E_\gamma,
\]

rather than the heat equation $\partial_s\gamma=D_t\dot\gamma$.
The former is an ODE on the $H^1(J,N)$ manifold; the latter is a parabolic
PDE requiring a different existence and regularity theory. Our construction
also does not imply that an $H^1(J,N)$ initial curve becomes instantly smooth:
smoothness will be proved for the critical points obtained as limits.

\end{remark}

\begin{theorem}[Weak critical points and smooth geodesics]

\label{teo:criticos-h1-geodesicas-suaves}

A curve $\gamma\in\mathcal P$ is a critical point of $\mathcal E$
if and only if it is a smooth constant-speed geodesic, with fixed
endpoints or closed, according to the choice of $\mathcal P$.
In the closed case, velocities also agree across the endpoint identification.

\end{theorem}

\begin{proof}

Suppose $d\mathcal E_\gamma=0$. On a coordinate segment $J_0\subseteq J$, write
$x=(x^1,\ldots,x^n)$ for $\gamma$ and use test fields supported in
the segment's interior. The coordinate first-variation formula implies,
for every $b\in\{1,\ldots,n\}$,

\begin{equation}\label{eq:regularidad-geodesica-momento}
 p_b:=\displaystyle\sum_{a=1}^n h_{ab}(x)(x^a)'\in L^2(J_0),
 \qquad
 p_b'=\frac12\displaystyle\sum_{a,c=1}^n
                    \partial_bh_{ac}(x)(x^a)'(x^c)'\in L^1(J_0).
\end{equation}

The last assertion uses $x'\in L^2(J_0,\mathbb R^n)$ and bounded coefficients over the
segment's compact image. An $L^1(J_0)$ function whose distributional derivative
belongs to $L^1(J_0)$ has an absolutely continuous representative: subtracting
the integral of its derivative leaves a distribution with zero derivative,
hence a constant. Thus $p_b$ is continuous and bounded on smaller closed
segments. Since

\[
 (x^a)'=\displaystyle\sum_{b=1}^n h^{ab}(x)p_b,
\]

the right-hand side is continuous and $x$ is $C^1$. The right-hand side
of \eqref{eq:regularidad-geodesica-momento} is now continuous, so $p$
is $C^1$ and $x$ is $C^2$. The two equations form a smooth first-order
system for $(x,p)$; induction gives $x,p\in C^\infty$.
The same argument reaches the endpoints of $I$: the terms in
\eqref{eq:regularidad-geodesica-momento} are integrable up to them,
so their primitives have limits; induction then applies to one-sided
derivatives. The smooth system with those endpoint data extends the
curve to a neighborhood of each endpoint. On $S^1$, local regularity
in charts crossing the join includes agreement of all derivatives.

We may now integrate by parts in
\eqref{eq:primera-variacion-curvas-h1}. Interior-supported fields give
$D_t\dot\gamma=0$. Metric compatibility yields
$\frac d{dt}|\dot\gamma|_{\mathbf h}^2=0$. Conversely, integration
by parts along a smooth geodesic gives
$d\mathcal E_\gamma(\mathbf V)=0$: the interior term vanishes, as do
the boundary terms by fixed endpoints or periodicity. The identity
extends to $H^1(J,\gamma^*TN)$ fields by continuity and density.

\end{proof}

\subsection{The Palais--Smale condition for energy}
\label{subsec:ps-energia-curvas}

\begin{theorem}[Compactness of almost critical sequences]

\label{teo:palais-smale-energia-curvas}

Energy $\mathcal E$ satisfies Palais--Smale with respect to $\mathbf G$
on $\mathcal P_{p,q}(N)$ when $N$ is complete, and on $\Lambda(N)$
when $N$ is compact.

\end{theorem}

\begin{proof}

Let $(\gamma_n)$ satisfy $\mathcal E(\gamma_n)\leq a$ and
$\|d\mathcal E_{\gamma_n}\|_{\mathbf G'}\to0$. For $t,\tau$ on one
parametrized arc, Cauchy--Schwarz gives

\[
 d_{\mathbf h}(\gamma_n(t),\gamma_n(\tau))
 \leq \sqrt{2a}\,|t-\tau|^{1/2}.
\]

For fixed endpoints, all images lie in the closed ball centered at $p$
of radius $\sqrt{2a}$, which is compact by
Theorem~\ref{teo:variedades-riemannianas-hopf-rinow}. For loops, use
compactness of $N$ directly. Arzelà--Ascoli gives a subsequence converging
uniformly to a continuous curve $\gamma$. The weak compactness argument
on coordinate segments used in the proof of
Theorem~\ref{teo:completitud-metrica-h1-curvas} shows that
$\gamma\in H^1(J,N)$.

In an exponential chart centered at $\gamma$, put
$\mathbf u_n=\Phi_\gamma(\gamma_n)$. The weak chain rule for the normal
logarithm bounds $\|\mathbf u_n\|_{H^1(J,\gamma^*TN)}$ by a constant times
$1+\|\dot\gamma\|_{L^2(J,\gamma^*TN)}+\|\dot\gamma_n\|_{L^2(J,\gamma_n^*TN)}$.
Thus $(\mathbf u_n)$ is bounded in $H^1(J,\gamma^*TN)$ and
$\|\mathbf u_n\|_{C^0}\to0$. After further extraction, its weak $H^1(J,\gamma^*TN)$
limit is zero and $D_t\mathbf u_n\rightharpoonup0$ in $L^2(J,\gamma^*TN)$.
The weak limit is identified using the already known uniform limit.
For fixed endpoints, $\mathbf u_n$ vanishes at both endpoints.

The decisive step is to turn weak convergence into strong convergence.
In \eqref{eq:velocidad-exp-curvas}, abbreviate
$A_n=A(\gamma,\mathbf u_n)$ and $B_n=B(\gamma,\mathbf u_n)$.
Differentiating $\mathcal E\circ R_\gamma$ at $\mathbf u_n$ in the
direction $\mathbf u_n$ gives exactly

\begin{align}
 D(\mathcal E\circ R_\gamma)(\mathbf u_n)[\mathbf u_n]
 &=\int_J|A_nD_t\mathbf u_n|_{\mathbf h}^2dt
   +\int_J\mathbf h(B_n\dot\gamma,A_nD_t\mathbf u_n)dt+\varepsilon_n,
 \label{eq:ps-energia-termino-coercivo}\\
 |\varepsilon_n|&\leq C\|\mathbf u_n\|_{C^0}
      \bigl(\|D_t\mathbf u_n\|_{L^2(J,\gamma^*TN)}^2+\|\dot\gamma\|_{L^2(J,\gamma^*TN)}^2\bigr).
 \label{eq:ps-energia-resto}
\end{align}

To check the remainder, the field $w_n=A_nD_t\mathbf u_n+B_n\dot\gamma$
is paired with
$(D_vA_n[\mathbf u_n])D_t\mathbf u_n+
(D_vB_n[\mathbf u_n])\dot\gamma$.
The vertical derivatives of $A$ and $B$ are bounded. Cauchy--Schwarz and
$2ab\leq a^2+b^2$ give \eqref{eq:ps-energia-resto}, so
$\varepsilon_n\to0$. Also, $A_n,B_n\to\operatorname{id}$ uniformly. Set

\[
\begin{aligned}
 \alpha_n&=\|A_n-\operatorname{id}\|_{L^\infty(J)},\qquad
 \beta_n=\|B_n-\operatorname{id}\|_{L^\infty(J)},\\
 L&=\sup_{n\in\mathbb N}\|D_t\mathbf u_n\|_{L^2(J,\gamma^*TN)}<\infty.
\end{aligned}
\]
Cauchy--Schwarz gives
\[
\begin{aligned}
 &\left|\int_J\mathbf h(B_n\dot\gamma,A_nD_t\mathbf u_n)\,dt
       -\int_J\mathbf h(\dot\gamma,D_t\mathbf u_n)\,dt\right|\\
 &\qquad\leq[\beta_n(1+\alpha_n)+\alpha_n]
                 \|\dot\gamma\|_{L^2(J,\gamma^*TN)}L\longrightarrow0.
\end{aligned}
\]
Since $\dot\gamma\in L^2(J,\gamma^*TN)$ is fixed and
$D_t\mathbf u_n\rightharpoonup0$, the second integral converges to zero.
The preceding bound proves the same for the first one.

The left-hand side of \eqref{eq:ps-energia-termino-coercivo} also tends
to zero. Indeed, put
$\mathbf W_n=dR_\gamma(\mathbf u_n)\mathbf u_n$.
On the chosen normal tube, the vertical component of $d\exp$ and its
first covariant derivatives are bounded. Differentiating $\mathbf W_n$
along $\gamma_n$, the chain rule gives, for almost every $t\in J$,

\[
 \begin{aligned}
 |\mathbf W_n|&\leq C|\mathbf u_n|,\\
 |D_t\mathbf W_n|&\leq C\bigl(
 |D_t\mathbf u_n|+
 (|\dot\gamma|+|D_t\mathbf u_n|)|\mathbf u_n|\bigr).
 \end{aligned}
\]

The constant $C$ depends on the tube over the compact image of
$\gamma$, but not on $n$. Pointwise norms are taken in the
corresponding fiber. Integrating and increasing $C$ if necessary yields

\[
 \begin{aligned}
 \|\mathbf W_n\|_{\mathbf G_{\gamma_n}}
 \leq C\bigl[&
 (1+\|\mathbf u_n\|_{C^0(J)})\|D_t\mathbf u_n\|_{L^2(J,\gamma^*TN)}\\
 &+(1+\|\dot\gamma\|_{L^2(J,\gamma^*TN)})\|\mathbf u_n\|_{C^0(J)}\bigr].
 \end{aligned}
\]

The bounds already obtained make this right-hand side uniformly bounded;
smallness of $\mathbf u_n$ in $H^1(J,\gamma^*TN)$ is not required.
By the Palais--Smale hypothesis,

\[
 |D(\mathcal E\circ R_\gamma)(\mathbf u_n)[\mathbf u_n]|
 =|d\mathcal E_{\gamma_n}(\mathbf W_n)|
 \leq\|d\mathcal E_{\gamma_n}\|_{\mathbf G_{\gamma_n}'}
       \|\mathbf W_n\|_{\mathbf G_{\gamma_n}}
 \longrightarrow0.
\]

The uniformly bounded inverses of $A_n$ in
\eqref{eq:ps-energia-termino-coercivo} now imply
$\|D_t\mathbf u_n\|_{L^2(J,\gamma^*TN)}\to0$.
Together with $\|\mathbf u_n\|_{C^0}\to0$, this gives
$\mathbf u_n\to0$ in $H^1(J,\gamma^*TN)$, hence $\gamma_n\to\gamma$ in
$\mathcal P$. Continuity of the differential makes $\gamma$ critical,
and Theorem~\ref{teo:criticos-h1-geodesicas-suaves} makes it smooth.

\end{proof}

\begin{remark}[Scope of the compactness statement]

\label{obs:limites-ps-curvas-potencial}

Completeness of $N$ does not suffice for free loops: in $N=\mathbb R$,
the constant loops $\gamma_n(t)=n$ have zero energy and gradient but
no convergent subsequence. Nor can an arbitrary smooth potential be
added while retaining the fixed-endpoint statement on a complete
noncompact manifold without further hypotheses. On
$\mathcal P_{0,0}(\mathbb R)$, the functional

\[
 \mathcal A(\gamma)=\frac12\int_0^1
           \bigl(|\gamma'(t)|^2-\pi^2|\gamma(t)|^2\bigr)dt
\]

has the curves $\gamma_n(t)=n\sin(\pi t)$ as critical points of value
zero. Indeed, $\gamma_n''+\pi^2\gamma_n=0$ and the endpoints vanish,
so integration by parts gives $d\mathcal A_{\gamma_n}=0$.
Their $H^1(J,\mathbb R^n)$ norms tend to infinity. Coercivity hypotheses or another
mechanism controlling excursions and kinetic energy are needed;
smoothness of the potential does not provide this control.

\end{remark}

\begin{theorem}[Obtaining geodesics by continuous descent]

\label{teo:geodesicas-descenso-continuo}

Let $\gamma_0\in\mathcal P_{p,q}(N)$ with $N$ complete, or
$\gamma_0\in\Lambda(N)$ with $N$ compact. The problem
$\partial_s\gamma=-\operatorname{grad}_{\mathbf G}\mathcal E(\gamma)$,
$\gamma(0)=\gamma_0$, has a unique solution for all $s\geq0$.
Energy decreases, and there are $s_n\to\infty$ such that
$\gamma(s_n)$ converges in $H^1(J,N)$ to a smooth geodesic $\gamma_*$
in the same component. If \eqref{eq:desigualdad-gradiente-local} holds
near $\gamma_*$, the whole trajectory converges to $\gamma_*$.
In particular, this holds if its Hessian is an isomorphism between the
tangent space and its dual.

\end{theorem}

\begin{proof}

Theorems~\ref{teo:completitud-metrica-h1-curvas} and
\ref{teo:palais-smale-energia-curvas} verify exactly the hypotheses of
Theorem~\ref{teo:descenso-global-palais-smale}. Then apply
Theorem~\ref{teo:criticos-h1-geodesicas-suaves}.
The limit lies in the same component: for large $n$, $\gamma(s_n)$
and $\gamma_*$ lie in one exponential chart and the chart segment
joins them. The descent trajectory joins $\gamma_0$ to $\gamma(s_n)$.
The full convergence criterion and its corollary give the last assertions.

\end{proof}

This result finds a critical point, not necessarily a global energy
minimum. Starting exactly at an unstable geodesic produces a constant
trajectory at that critical point. Nor does it guarantee obtaining all
geodesics from a single initial curve. In the loop space, a nontrivial
free homotopy class excludes a constant limit: the trajectory and chart
homotopies are continuous after evaluation at $t$. In the trivial class,
the limit may be a constant loop. Forcing other critical points requires
topological constraints or minimax procedures, not merely decreasing
energy.

\subsection{Hessian, Jacobi fields, and convergence}
\label{subsec:hessiana-descenso-geodesicas}

\begin{proposition}[The Hessian as a compact perturbation of the identity]

\label{prop:hessiana-geodesicas-fredholm}

At a geodesic $\gamma\in\mathcal P$, the Hessian of $\mathcal E$ is

\begin{equation}\label{eq:hessiana-curvas-geodesicas}
 d^2\mathcal E_\gamma(\mathbf V,\mathbf W)
 =\int_J\left(\mathbf h(D_t\mathbf V,D_t\mathbf W)
 -\mathbf h(R(\mathbf V,\dot\gamma)\dot\gamma,\mathbf W)\right)dt.
\end{equation}

The bounded self-adjoint operator $A_\gamma$ representing it through
$\mathbf G_\gamma$ has the form $\operatorname{id}-K_\gamma$, where
$K_\gamma$ is compact. Its kernel consists of smooth Jacobi fields with
the corresponding boundary conditions. In particular, $A_\gamma$ is
invertible if and only if this kernel is zero.

\end{proposition}

\begin{proof}

For smooth fields, differentiate
$\displaystyle \int_J\mathbf h(\partial_t\alpha,D_t\mathbf W)dt$ along a variation
$\alpha$ whose variation velocity is $\mathbf V$.
Torsion-freeness gives $D_s\partial_t\alpha=D_t\partial_s\alpha$, and
the curvature identity gives
$D_sD_t\mathbf W=D_tD_s\mathbf W+R(\partial_s\alpha,\partial_t\alpha)\mathbf W$.
The term involving $D_s\mathbf W$ vanishes after integration by parts:
$D_t\dot\gamma=0$ and the endpoint terms vanish by the imposed conditions.
Curvature symmetries give the term in
\eqref{eq:hessiana-curvas-geodesicas}. Both sides are continuous bilinear
forms on $H^1(J,\gamma^*TN)$, so density of smooth fields with the same boundary
conditions extends the identity to all tangent vectors.

Define $K_\gamma$ by Riesz representation through

\[
 \mathbf G_\gamma(K_\gamma\mathbf V,\mathbf W)
 =\int_J\mathbf h\bigl(\mathbf V+R(\mathbf V,\dot\gamma)\dot\gamma,
                       \mathbf W\bigr)dt.
\]

The inclusion $H^1(J,\gamma^*TN)\hookrightarrow L^2(J,\gamma^*TN)$
is compact by Rellich--Kondrashov. The endomorphism
$\mathbf V\mapsto\mathbf V+R(\mathbf V,\dot\gamma)\dot\gamma$
is bounded on $L^2(J,\gamma^*TN)$, since $\gamma$ is smooth and $J$ compact.
The inclusion $L^2(J,\gamma^*TN)\to\mathcal H_\gamma'$ is continuous.
Thus $K_\gamma$ is a composition of a compact operator with continuous
operators. The Hessian formula gives
$A_\gamma=\operatorname{id}-K_\gamma$.

If $A_\gamma\mathbf V=0$, interior test fields give
$D_t^2\mathbf V+R(\mathbf V,\dot\gamma)\dot\gamma=0$ weakly.
The one-dimensional regularity argument used for
\eqref{eq:gradiente-geodesicas-green} makes $\mathbf V$ smooth.
The boundary conditions can be made explicit. For fixed endpoints,
$\mathbf V(0)=\mathbf V(1)=0$. For loops, the weak equation is
interpreted on $S^1$ and regularity also applies in a chart crossing the
identified endpoints; hence both $\mathbf V$ and $D_t\mathbf V$
are periodic. Equivalently, after regularity on the interval has been
established, periodic test fields force the boundary term to satisfy
$D_t\mathbf V(1)=D_t\mathbf V(0)$. The converse follows by integration
by parts with these conditions. Suppose $\ker A_\gamma=\{0\}$. If no $c>0$ satisfied
$\|A_\gamma\mathbf V\|\geq c\|\mathbf V\|$, there would be unit vectors
$\mathbf V_n$ with $A_\gamma\mathbf V_n\to0$.
A subsequence of $K_\gamma\mathbf V_n$ converges; since
$\mathbf V_n=A_\gamma\mathbf V_n+K_\gamma\mathbf V_n$, it converges to
a unit vector in the kernel, a contradiction. Thus the image is closed.
By self-adjointness, its orthogonal complement is the kernel, so the
image is the whole space. The lower bound gives a bounded inverse.

\end{proof}

For fixed endpoints, absence of nonzero Jacobi fields vanishing at both
endpoints means that the endpoints are not conjugate along $\gamma$.
If a descent subsequence converges to such a geodesic,
Corollary~\ref{cor:descenso-limite-no-degenerado} gives exponential
convergence of the entire trajectory. For a nonconstant closed geodesic,
$\dot\gamma$ is a nonzero periodic Jacobi field, corresponding to a
shift of the parameter origin. Thus the corollary cannot be applied
automatically in the loop space; full convergence requires, for example,
verification of \eqref{eq:desigualdad-gradiente-local} or additional
analysis of the critical family.

\subsection{Discrete descent with step search}
\label{subsec:descenso-discreto-geodesicas}

\begin{theorem}[Descent iteration with the Armijo condition]

\label{teo:descenso-discreto-armijo-geodesicas}

Under the hypotheses of Theorem~\ref{teo:geodesicas-descenso-continuo},
fix an initial curve $\gamma_0$, $a=\mathcal E(\gamma_0)$, and
$0<\sigma,\beta<1$. If $a=0$, $\gamma_0$ is already critical.
If $a>0$, choose $\rho_a$ as in
Lemma~\ref{lem:exp-curvas-estimaciones-uniformes}, for the compact set
containing all curves of energy at most $a$, and fix

\[
 0<\alpha_{\max}\leq\frac{\rho_a}{\sqrt{2a}}.
\]

Given $\gamma_n$, compute $\mathbf Z_n$ using
\eqref{eq:gradiente-geodesicas-forma-debil}. If $\mathbf Z_n=0$, stop.
Otherwise, choose the first integer $k_n\geq0$ such that, with
$\alpha_n=\alpha_{\max}\beta^{k_n}$,

\begin{equation}\label{eq:armijo-energia-curvas}
 \mathcal E\bigl(R_{\gamma_n}(-\alpha_n\mathbf Z_n)\bigr)
 \leq\mathcal E(\gamma_n)-\sigma\alpha_n\|\mathbf Z_n\|_{\mathbf G_{\gamma_n}}^2,
 \qquad
 \gamma_{n+1}=R_{\gamma_n}(-\alpha_n\mathbf Z_n).
\end{equation}

The search terminates at every iteration. If the procedure does not
stop after finitely many steps, the energies converge,
$\|\mathbf Z_n\|_{\mathbf G_{\gamma_n}}\to0$, and a subsequence converges
in $H^1(J,N)$ to a smooth geodesic in the same component. Accepted steps
are bounded below by one positive constant independent of $n$.

\end{theorem}

\begin{proof}

The compact set in the statement is $N$ for loops. For fixed endpoints,
it is the closed ball centered at $p$ of radius $\sqrt{2a}$.
The lemma applies over this compact set and a normal tube around it.
Inductively, the Armijo condition keeps $\mathcal E(\gamma_n)\leq a$.
Bound \eqref{eq:cota-gradiente-energia-curvas} and the choice of
$\alpha_{\max}$ ensure that every trial step stays in the retraction
domain and in the ball where
\eqref{eq:taylor-uniforme-energia-retraccion} holds.

Taylor's formula with integral remainder gives, for
$0<\alpha\leq\alpha_{\max}$,

\[
 \mathcal E\bigl(R_{\gamma_n}(-\alpha\mathbf Z_n)\bigr)
 \leq\mathcal E(\gamma_n)-\alpha\|\mathbf Z_n\|_{\mathbf G}^2
       +\frac{L_a}{2}\alpha^2\|\mathbf Z_n\|_{\mathbf G}^2.
\]

Here we used $dR_{\gamma_n}(0)=\operatorname{id}$ and the definition
of gradient. Taking $L_a>0$, every
$\alpha\leq2(1-\sigma)/L_a$ satisfies Armijo. Since trial steps form
a geometric sequence tending to zero, the search terminates.
If $\alpha_{\max}$ is not accepted, the first step crossing this
threshold exceeds $\beta\,2(1-\sigma)/L_a$. In all cases,

\[
 \alpha_n\geq\underline\alpha:=
 \min\left\{\alpha_{\max},\frac{2\beta(1-\sigma)}{L_a}\right\}>0.
\]

Summing \eqref{eq:armijo-energia-curvas} from $n=0$ to $m$ gives

\[
 \sigma\underline\alpha\displaystyle\sum_{n=0}^{m}
       \|\mathbf Z_n\|_{\mathbf G_{\gamma_n}}^2
 \leq\mathcal E(\gamma_0)-\mathcal E(\gamma_{m+1})\leq a.
\]

Thus the series of squared norms converges and the norms tend to zero.
Palais--Smale gives a convergent subsequence, and critical-point regularity
gives a smooth geodesic. Each step is joined to the preceding curve by
$u\mapsto R_{\gamma_n}(-u\alpha_n\mathbf Z_n)$, $0\leq u\leq1$.
This preserves endpoints and the component; the limit belongs to it
by the chart argument used earlier.

\end{proof}

\begin{proposition}[Convergence of the whole discrete sequence]

\label{prop:convergencia-total-armijo}

If a subsequence of the preceding iteration converges to $\gamma_*$
and \eqref{eq:desigualdad-gradiente-local} holds near $\gamma_*$, then
$\gamma_n\to\gamma_*$ in $H^1(J,N)$. If $\theta=1/2$, energies converge
geometrically, as do distances, possibly with a different ratio.

\end{proposition}

\begin{proof}

If the procedure terminates, extend the sequence constantly from its
last term; the conclusion is immediate. Suppose it does not terminate.
Let $e_n=\mathcal E(\gamma_n)-\mathcal E(\gamma_*)$.
Monotonicity of energy and convergence of the subsequence give
$e_n\geq0$ and $\displaystyle\lim_{n\to\infty}e_n=0$.
If some $e_n=0$, Armijo would force $\mathbf Z_n=0$, so in the case
under consideration $e_n>0$ for every $n$.
The length of the retraction segment and
\eqref{eq:diferencial-retraccion-curvas} give

\[
 d_{\mathbf G}(\gamma_{n+1},\gamma_n)
 \leq C_a\alpha_n\|\mathbf Z_n\|_{\mathbf G_{\gamma_n}}.
\]

Let $U$ be a neighborhood of $\gamma_*$ where
\eqref{eq:desigualdad-gradiente-local} holds, with constants $C>0$ and
$0<\theta\leq1/2$. As long as $\gamma_n\in U$, concavity of
$t^\theta$ and Armijo imply

\[
 \begin{aligned}
 e_n^\theta-e_{n+1}^\theta
 &\geq\theta e_n^{\theta-1}(e_n-e_{n+1})\\
 &\geq\frac{\theta\sigma}{C}
             \alpha_n\|\mathbf Z_n\|_{\mathbf G_{\gamma_n}}\\
 &\geq\frac{\theta\sigma}{CC_a}
                   d_{\mathbf G}(\gamma_{n+1},\gamma_n).
 \end{aligned}
\]

Put $S=CC_a/(\theta\sigma)$ and choose $\rho>0$ such that
$B_{d_{\mathbf G}}(\gamma_*,\rho)\subset U$. This is possible because
the metric is strong and its distance induces the manifold topology.
From the convergent subsequence choose an index $n_0$ such that

\[
 d_{\mathbf G}(\gamma_{n_0},\gamma_*)<\frac\rho3,
 \qquad S e_{n_0}^{\theta}<\frac\rho3.
\]

If $m>n_0$ were the first index with
$\gamma_m\notin B_{d_{\mathbf G}}(\gamma_*,\rho)$, all terms with
indices $n_0,\ldots,m-1$ would lie in $U$. Summing the preceding
estimate, including the last step, would give

\[
 \begin{aligned}
 d_{\mathbf G}(\gamma_m,\gamma_*)
 &\leq d_{\mathbf G}(\gamma_{n_0},\gamma_*)
       +\sum_{j=n_0}^{m-1}d_{\mathbf G}(\gamma_{j+1},\gamma_j)\\
 &\leq d_{\mathbf G}(\gamma_{n_0},\gamma_*)
       +S(e_{n_0}^{\theta}-e_m^{\theta})
 <\frac{2\rho}{3},
 \end{aligned}
\]

contradicting the choice of $m$. Thus the entire tail remains in that
ball. For $m>n\geq n_0$, the same argument yields

\[
 \sum_{j=n}^{m-1}d_{\mathbf G}(\gamma_{j+1},\gamma_j)
 \leq S(e_n^\theta-e_m^\theta)\leq S e_n^\theta.
\]

Since $e_n\to0$, the sequence is Cauchy. Completeness of the component
gives a limit, which equals $\gamma_*$ by the convergent subsequence.
Passing to the limit as $m\to\infty$ gives

\begin{equation}\label{eq:cola-distancias-armijo}
 d_{\mathbf G}(\gamma_n,\gamma_*)
 \leq\sum_{j=n}^{\infty}d_{\mathbf G}(\gamma_{j+1},\gamma_j)
 \leq\frac{CC_a}{\theta\sigma}\,e_n^\theta,
 \qquad n\geq n_0.
\end{equation}

If $\theta=1/2$, the gradient inequality gives
$\|\mathbf Z_n\|_{\mathbf G_{\gamma_n}}^2\geq e_n/C^2$.
Increase $C$, if necessary, before choosing $n_0$, so that
$q:=1-\sigma\underline\alpha/C^2$ belongs to $(0,1)$.
Armijo and the lower bound on step sizes give $e_{n+1}\leq q e_n$.
By induction and \eqref{eq:cola-distancias-armijo}, for
$n\geq n_0$,

\[
 e_n\leq e_{n_0}q^{n-n_0},\qquad
 d_{\mathbf G}(\gamma_n,\gamma_*)
 \leq\frac{2CC_a}{\sigma}\sqrt{e_{n_0}}\,q^{(n-n_0)/2}.
\]

This proves both convergence rates.\qedhere

\end{proof}

\begin{example}[The Euclidean case]

\label{ej:descenso-geodesicas-euclidiano}

For $N=\mathbb R^d$, write
$\gamma(t)=p+t(q-p)+v(t)$, where $v\in H^1_0(I,\mathbb R^d)$.
The segment $\gamma_*(t)=p+t(q-p)$ is the only critical point.
Indeed,

\[
 \mathcal E(\gamma)=\frac12|q-p|^2+\frac12\int_0^1|v'(t)|^2dt,
\]

because the cross term is $(q-p)\cdot\displaystyle \int_0^1v'(t)dt=0$.
In the sine expansion $v(t)=\displaystyle \sum_{k=1}^\infty v_k\sin(k\pi t)$,
the weak gradient equation determines

\[
 Z_k=\lambda_kv_k,\qquad
 \lambda_k=\frac{k^2\pi^2}{1+k^2\pi^2},\qquad
 v_k(s)=e^{-\lambda_ks}v_k(0).
\]

To verify this, test \eqref{eq:gradiente-geodesicas-forma-debil} against
each sine function: orthogonality diagonalizes both the zeroth-order
and first-order terms. Convergence of the series in $H^1_0((0,1),\mathbb R^d)$ and
$0<\pi^2/(1+\pi^2)\leq\lambda_k<1$ allow passage from finite sums
to the limit. Thus

\[
 \|\gamma(s)-\gamma_*\|_{H^1((0,1),\mathbb R^d)}
 \leq e^{-\pi^2s/(1+\pi^2)}\|\gamma_0-\gamma_*\|_{H^1((0,1),\mathbb R^d)}.
\]

Notice that high frequencies decay at rates tending to one, rather
than at rates $k^2\pi^2$ as for the heat equation. Sobolev descent does
not have the instantaneous regularizing effect of that parabolic flow.
The discrete iteration acts on the same coefficient by
$v_{k,n+1}=(1-\alpha_n\lambda_k)v_{k,n}$; Armijo search controls energy
decrease simultaneously at all frequencies.

\end{example}

In a numerical implementation,
\eqref{eq:gradiente-geodesicas-forma-debil} is approximated using a
finite basis or finite elements on the interval. The preceding theorems
concern the exact $H^1(J,N)$ gradient and exact energy evaluations.
Convergence when the gradient, exponential map, or energy is computed
inexactly requires error estimates and does not follow automatically
from Armijo.

\section{Yang--Mills theory and gauge reduction}
\label{sec:aplicacion-yang-mills}
\index{Yang--Mills}
\index{gauge group}
\index{Coulomb condition}

The second application has an additional symmetry. A connection is not an
isolated object: the gauge group transforms distinct connections describing
the same geometry. At critical connections, this invariance places the directions tangent
to gauge orbits in the kernel of the Hessian of the Yang--Mills
functional and leads naturally to the Coulomb condition.

Let \((M,\mathbf{g})\) be a compact Riemannian manifold of dimension $m$ with or without boundary,
and let \(\mathbf{P}\longrightarrow M\) be a principal bundle with compact structure
group \(G\). Denote its adjoint bundle by
\[
 \operatorname{ad}\mathbf{P}
 :=
 \mathbf{P}\times_{\operatorname{Ad}}\mathfrak g
\]
and fix an \(\operatorname{Ad}\)-invariant inner product on \(\mathfrak g\).
We denote the induced metric on \(\operatorname{ad}\mathbf{P}\) by \(\mathbf{h}_{\mathrm{ad}}\).

The curvature of a connection \(\mathbf{A}\) is an adjoint-bundle-valued \(2\)-form:
\[
 \mathbf{F}_{\mathbf{A}}\in\Omega^2(M,\operatorname{ad}\mathbf{P})
 :=
 \Gamma\left(
 \Lambda^2(T^*M)\otimes\operatorname{ad}\mathbf{P}
 \right).
\]
The superscript \(2\) indicates the degree of the form: for each \(x\in M\),
\(\mathbf{F}_{\mathbf{A}}(x)\) takes two vectors in \(T_xM\), depends
alternatingly on them, and takes values in \((\operatorname{ad}\mathbf{P})_x\).

The space of connections is affine, modeled on
\(\Omega^1(M,\operatorname{ad}\mathbf{P})\). If \(\mathbf{a}\) is an
\(\operatorname{ad}\mathbf{P}\)-valued \(1\)-form, the curvature of
\(\mathbf{A}+\mathbf{a}\) satisfies
\begin{equation}
\label{eq:curvatura-variada-ym}
 \mathbf{F}_{\mathbf{A}+\mathbf{a}}
 =
 \mathbf{F}_{\mathbf{A}}+d_{\mathbf{A}}\mathbf{a}
 +\frac{1}{2}[\mathbf{a}\wedge \mathbf{a}].
\end{equation}
The Yang--Mills functional is
\begin{equation}
\label{eq:funcional-yang-mills}
 \operatorname{YM}(\mathbf{A})
 :=
 \frac{1}{2}\int_M|\mathbf{F}_{\mathbf{A}}|_{\mathbf{g},\mathbf{h}_{\mathrm{ad}}}^2\,d\lambda_{\mathbf{g}}.
\end{equation}

\begin{proposition}[First variation of Yang--Mills]
\label{prop:primera-variacion-yang-mills}
Let $(M,\mathbf{g})$ be a compact Riemannian manifold with or without boundary,
let $\mathbf{P}\to M$ be a principal bundle with compact structure group, and let
$\mathbf{A}$ be a smooth connection. For
\(\mathbf{a}\in\Omega^1(M,\operatorname{ad}\mathbf{P})\),
\begin{equation}
\label{eq:primera-variacion-yang-mills}
 D\operatorname{YM}(\mathbf{A})[\mathbf{a}]
 =
 \int_M\langle d_{\mathbf{A}}^*\mathbf{F}_{\mathbf{A}},\mathbf{a}\rangle_{\mathbf{g},\mathbf{h}_{\mathrm{ad}}}\,d\lambda_{\mathbf{g}}
 +\int_{\partial M}
 \langle\iota_{\boldsymbol{\nu}}\mathbf{F}_{\mathbf{A}},\mathbf{a}\rangle_{\mathbf{g},\mathbf{h}_{\mathrm{ad}}}\,d\sigma_{\mathbf{g}}.
\end{equation}
Here $\boldsymbol{\nu}$ is the outward unit normal and
\[
 d_{\mathbf{A}}^*\mathbf{F}
 =-\operatorname{tr}_{\mathbf{g}}(\nabla^{\mathbf{A}}\mathbf{F})
\]
on \(\operatorname{ad}\mathbf{P}\)-valued \(2\)-forms.
Contraction is performed between the derivative index and the first index
of the form: for $\mathbf X\in T_xM$,
\[
 (d_{\mathbf A}^*\mathbf F)(\mathbf X)
 =-\sum_{i=1}^m(\nabla^{\mathbf A}_{\mathbf e_i}\mathbf F)(\mathbf e_i,\mathbf X).
\]
\end{proposition}

\begin{proof}
Substituting \(t\mathbf{a}\) into \eqref{eq:curvatura-variada-ym}, the part linear in
\(t\) is \(d_{\mathbf{A}}\mathbf{a}\). Consequently,
\[
 \left.\frac{d}{dt}\right|_{t=0}\operatorname{YM}(\mathbf{A}+t\mathbf{a})
 =
 \int_M\langle \mathbf{F}_{\mathbf{A}},d_{\mathbf{A}}\mathbf{a}\rangle_{\mathbf{g},\mathbf{h}_{\mathrm{ad}}}\,d\lambda_{\mathbf{g}}.
\]
The Green form from Chapter~\ref{cap:sobolev-haces} gives the adjoint and
the boundary term. To make the contraction explicit, define
\[
 \alpha(\mathbf X):=\sum_{j=1}^m
 \langle\mathbf F_{\mathbf A}(\mathbf X,\mathbf e_j),\mathbf a(\mathbf e_j)\rangle_{\mathbf h_{\mathrm{ad}}}.
\]
In an orthonormal frame normal at the point of calculation, the Leibniz rule and
antisymmetry of $\mathbf F_{\mathbf A}$ give
$\operatorname{div}\alpha^\sharp
=\langle\mathbf F_{\mathbf A},d_{\mathbf A}\mathbf a\rangle
 -\langle d_{\mathbf A}^*\mathbf F_{\mathbf A},\mathbf a\rangle$.
Using the norm of $2$-forms means that we sum once over
$i<j$, rather than twice over ordered pairs. Upon integration,
$\alpha(\boldsymbol\nu)=\langle\iota_{\boldsymbol\nu}\mathbf F_{\mathbf A},\mathbf a\rangle$
yields exactly \eqref{eq:primera-variacion-yang-mills}.
\end{proof}

If $M$ has no boundary, or if the induced tangential connection
on $\partial M$ is fixed, the admissible variations make the boundary term vanish.
In the second case, we require $\iota^*\mathbf a=0$, where
$\iota:\partial M\hookrightarrow M$ is the inclusion; the normal
component of $\mathbf a$ does not contribute because
$(\iota_{\boldsymbol\nu}\mathbf F_{\mathbf A})(\boldsymbol\nu)=0$.
The formal $L^2(M,T^*M\otimes\operatorname{ad}\mathbf P)$ representative of the differential and the interior
Euler--Lagrange equation are then
\begin{equation}
\label{eq:gradiente-ecuacion-yang-mills}
 \operatorname{grad}_{L^2(M,T^*M\otimes\operatorname{ad}\mathbf P)}\operatorname{YM}(\mathbf{A})
 =
 d_{\mathbf{A}}^*\mathbf{F}_{\mathbf{A}},
 \qquad
 d_{\mathbf{A}}^*\mathbf{F}_{\mathbf{A}}=0.
\end{equation}
If the connection on the boundary may also vary freely, the
critical point condition additionally includes
$\iota_{\boldsymbol\nu}\mathbf F_{\mathbf A}=0$ on $\partial M$.
Interior variations first determine the Yang--Mills equation;
arbitrariness of the tangential trace then determines this natural
boundary condition.

To write the Hessian, we first specify the algebraic term containing the
curvature. Define pointwise the endomorphism
\[
 \boldsymbol{\mathcal R}_{\mathbf{A}}\colon
 T^*M\otimes\operatorname{ad}\mathbf{P}
 \longrightarrow
 T^*M\otimes\operatorname{ad}\mathbf{P}
\]
by
\[
 \langle\boldsymbol{\mathcal R}_{\mathbf{A}}\eta,\zeta\rangle_{\mathbf{g},\mathbf{h}_{\mathrm{ad}};x}
 :=
 \langle \mathbf{F}_{\mathbf{A}},[\eta\wedge\zeta]\rangle_{\mathbf{g},\mathbf{h}_{\mathrm{ad}};x},
 \qquad
 \eta,\zeta\in T_x^*M\otimes(\operatorname{ad}\mathbf{P})_x.
\]
The Riesz representation theorem in each fiber determines a unique
\(\boldsymbol{\mathcal R}_{\mathbf{A}}\eta\). Its coefficients are
algebraic contractions of \(\mathbf{F}_{\mathbf{A}}\) with the bracket on
\(\mathfrak g\), so \(\boldsymbol{\mathcal R}_{\mathbf{A}}\) is a
smooth endomorphism of order zero. For $1\leq i,j\leq m$,
$[\eta\wedge\zeta](\mathbf e_i,\mathbf e_j)
=[\eta(\mathbf e_i),\zeta(\mathbf e_j)]
 -[\eta(\mathbf e_j),\zeta(\mathbf e_i)]$.
Antisymmetry of the bracket implies
$[\eta\wedge\zeta]=[\zeta\wedge\eta]$; thus the
endomorphism $\boldsymbol{\mathcal R}_{\mathbf A}$ is symmetric in each
fiber, as required for the Hessian.

\begin{proposition}[Hessian and gauge-fixed operator]
\label{prop:hessiana-simbolo-yang-mills}
Let $(M,\mathbf{g})$ be a compact Riemannian manifold with or without boundary,
let $\mathbf{P}\to M$ be a principal bundle with compact structure group, and let
$\mathbf{A}$ be a smooth connection. For
\(\mathbf{a},\mathbf{b}\in\Omega^1(M,\operatorname{ad}\mathbf{P})\),
\begin{equation}
\label{eq:segunda-variacion-bilineal-yang-mills}
 D^2\operatorname{YM}(\mathbf{A})[\mathbf{a},\mathbf{b}]
 =
 \int_M\left(
 \langle d_{\mathbf{A}}\mathbf{a},d_{\mathbf{A}}\mathbf{b}\rangle_{\mathbf{g},\mathbf{h}_{\mathrm{ad}}}
 +
 \langle \mathbf{F}_{\mathbf{A}},[\mathbf{a}\wedge\mathbf{b}]\rangle_{\mathbf{g},\mathbf{h}_{\mathrm{ad}}}
 \right)\,d\lambda_{\mathbf{g}}.
\end{equation}
If \(M\) is compact without boundary and \(\mathbf{A}\) is a Yang--Mills connection, the Hessian is
represented by
\begin{equation}
\label{eq:hessiana-yang-mills}
 H_{\mathbf{A}}\mathbf{a}
 =d_{\mathbf{A}}^*d_{\mathbf{A}}\mathbf{a}
 +\boldsymbol{\mathcal R}_{\mathbf{A}}\mathbf{a}.
\end{equation}
Moreover,
\[
 H_{\mathbf{A}}(d_{\mathbf{A}}\boldsymbol{\xi})=0
 \qquad
 \text{for every }\boldsymbol{\xi}\in\Omega^0(M,\operatorname{ad}\mathbf{P}).
\]
The gauge-fixed operator
\begin{equation}
\label{eq:operador-gauge-fijado-yang-mills}
 L_{\mathbf{A}}\mathbf{a}
 =
 (d_{\mathbf{A}}^*d_{\mathbf{A}}+d_{\mathbf{A}}d_{\mathbf{A}}^*)\mathbf{a}
 +\boldsymbol{\mathcal R}_{\mathbf{A}}\mathbf{a}
\end{equation}
has principal symbol
\begin{equation}
\label{eq:simbolo-principal-yang-mills-fijado}
 \boldsymbol{\sigma}_2(L_{\mathbf{A}})(x,\zeta)
 =
 -|\zeta|_{\mathbf{g}}^2
 \operatorname{id}_{T_x^*M\otimes(\operatorname{ad}\mathbf{P})_x},
 \qquad
 \zeta\in T_x^*M\setminus\{0\}.
\end{equation}
\end{proposition}

\begin{proof}
Substitute \(\mathbf{A}+t\mathbf{a}+s\mathbf{b}\) into
\eqref{eq:curvatura-variada-ym}. The coefficient
of \(ts\) in
\[
 \frac{1}{2}
 \|\mathbf{F}_{\mathbf{A}+t\mathbf{a}+s\mathbf{b}}\|_{
 L^2(M,\Lambda^2(T^*M)\otimes\operatorname{ad}\mathbf{P})}^2
\]
is
\[
 \int_M\langle d_{\mathbf{A}}\mathbf{a},d_{\mathbf{A}}\mathbf{b}\rangle_{\mathbf{g},\mathbf{h}_{\mathrm{ad}}}\,d\lambda_{\mathbf{g}}
 +
 \int_M\langle \mathbf{F}_{\mathbf{A}},[\mathbf{a}\wedge\mathbf{b}]\rangle_{\mathbf{g},\mathbf{h}_{\mathrm{ad}}}\,d\lambda_{\mathbf{g}},
\]
which proves \eqref{eq:segunda-variacion-bilineal-yang-mills}. When
\(M\) is compact without boundary, integration by parts in the first term gives
\eqref{eq:hessiana-yang-mills}.

Let \(\boldsymbol{\xi}\in\Omega^0(M,\operatorname{ad}\mathbf{P})\). Gauge
invariance implies, for every connection \(\mathbf{B}\),
\[
 D\operatorname{YM}(\mathbf{B})[d_{\mathbf{B}}\boldsymbol{\xi}]=0.
\]
Differentiate this identity at \(\mathbf{B}=\mathbf{A}\) in the direction
\(\mathbf{b}\). The term arising from the variation of
\(d_{\mathbf{B}}\boldsymbol{\xi}\) is multiplied by
\(D\operatorname{YM}(\mathbf{A})\) and vanishes because \(\mathbf{A}\) is critical. Thus
\[
 D^2\operatorname{YM}(\mathbf{A})[\mathbf{b},d_{\mathbf{A}}\boldsymbol{\xi}]=0
\]
for every \(\mathbf{b}\). Symmetry of the Hessian gives
\[
 H_{\mathbf{A}}(d_{\mathbf{A}}\boldsymbol{\xi})=0.
\]

Finally, for \(\zeta\in T_x^*M\setminus\{0\}\),
\[
 \boldsymbol{\sigma}_2(d_{\mathbf{A}}^*d_{\mathbf{A}})(x,\zeta)
 =
 -\iota_{\zeta^\sharp}(\zeta\wedge\cdot),
\]
whereas
\[
 \boldsymbol{\sigma}_2(d_{\mathbf{A}}d_{\mathbf{A}}^*)(x,\zeta)
 =
 -\zeta\wedge\iota_{\zeta^\sharp}(\cdot).
\]
The identity
\[
 \iota_{\zeta^\sharp}(\zeta\wedge a)
 +
 \zeta\wedge\iota_{\zeta^\sharp}a
 =
 |\zeta|_{\mathbf{g}}^2a
\]
shows that the sum of the two symbols is
\(-|\zeta|_{\mathbf{g}}^2\operatorname{id}\). Since
\(\boldsymbol{\mathcal R}_{\mathbf{A}}\) has order zero,
we obtain \eqref{eq:simbolo-principal-yang-mills-fijado}.
\end{proof}

On the compact manifold without boundary, the Coulomb condition
\[
 d_{\mathbf{A}}^*\mathbf{a}=0
\]
selects the \(L^2(M,T^*M\otimes\operatorname{ad}\mathbf P)\)-orthogonal complement of the infinitesimal directions
\(d_{\mathbf{A}}\boldsymbol{\xi}\). On this subspace, the term
\(d_{\mathbf{A}}d_{\mathbf{A}}^*\mathbf{a}\) vanishes, so
\(L_{\mathbf{A}}\) agrees with the Hessian. On the full space of \(1\)-forms,
this term completes the symbol and yields an elliptic operator.

On a compact base without boundary, projection onto the Coulomb
complement is obtained by solving an elliptic equation on adjoint-bundle-valued
\(0\)-forms.

By Proposition~\ref{prop:descomposicion-coulomb-compacta}, for every
$s\geq1$ the Coulomb projection is given by
\[
 \Pi_{\mathbf A}\mathbf a
 =\mathbf a-d_{\mathbf A}
 \left(\Delta_{\mathbf A}^{(0)}|_{(\ker d_{\mathbf A})^\perp}\right)^{-1}
 d_{\mathbf A}^*\mathbf a,
 \qquad
 \mathbf a\in H^s(M,T^*M\otimes\operatorname{ad}\mathbf P).
\]
The inverse acts from $H^{s-1}(M,\operatorname{ad}\mathbf P)\cap(\ker d_{\mathbf A})^\perp$ to
$H^{s+1}(M,\operatorname{ad}\mathbf P)\cap(\ker d_{\mathbf A})^\perp$, with orthogonality in
$L^2(M,\operatorname{ad}\mathbf P)$. Thus $\Pi_{\mathbf A}$ is continuous
on $H^s(M,T^*M\otimes\operatorname{ad}\mathbf P)$ and takes values in $\ker d_{\mathbf A}^*$.

Now fix
\[
 s>\frac{m}{2}+1
\]
and a smooth connection \(\mathbf{A}_0\). The Sobolev completion of the space of
connections is the affine space
\[
 \mathcal A^s(\mathbf{P})
 :=
 \mathbf{A}_0+
 H^s\left(
 M,T^*M\otimes\operatorname{ad}\mathbf{P}
 \right),
\]
and
\[
 T_{\mathbf{A}}\mathcal A^s(\mathbf{P})
 =
 H^s\left(
 M,T^*M\otimes\operatorname{ad}\mathbf{P}
 \right).
\]

Identity \eqref{eq:curvatura-variada-ym} shows that the curvature map
loses one derivative:
\[
 \mathcal F\colon
 \mathcal A^s(\mathbf{P})
 \longrightarrow
 H^{s-1}\left(
 M,\Lambda^2(T^*M)\otimes\operatorname{ad}\mathbf{P}
 \right),
 \qquad
 \mathcal F(\mathbf{A}):=\mathbf{F}_{\mathbf{A}}.
\]
The functional
\[
 \operatorname{YM}\colon
 \mathcal A^s(\mathbf{P})\longrightarrow\mathbb R,
 \qquad
 \operatorname{YM}(\mathbf{A})
 =
 \frac{1}{2}\int_M|\mathbf{F}_{\mathbf{A}}|_{\mathbf{g},\mathbf{h}_{\mathrm{ad}}}^2\,d\lambda_{\mathbf{g}},
\]
is also smooth. Indeed, the linear part of \(\mathcal F\) is a
first-order differential operator, and the remaining part is quadratic. The
Sobolev product estimates control
\([\mathbf{a}\wedge\mathbf{a}]\), while
the inclusion
\[
 H^{s-1}\left(
 M,\Lambda^2(T^*M)\otimes\operatorname{ad}\mathbf{P}
 \right)
 \hookrightarrow
 L^2\left(
 M,\Lambda^2(T^*M)\otimes\operatorname{ad}\mathbf{P}
 \right)
\]
allows us to integrate the squared norm of the curvature.

The gauge group $\mathcal G^{s+1}(\mathbf P)=H^{s+1}(M,\operatorname{Ad}(\mathbf P))$ acts smoothly on
\(\mathcal A^s(\mathbf{P})\). For each smooth connection \(\mathbf{A}\), the Hessian
and Coulomb operators induce continuous linear maps
\[
 H_{\mathbf{A}},L_{\mathbf{A}}\colon
 H^s\left(
 M,T^*M\otimes\operatorname{ad}\mathbf{P}
 \right)
 \longrightarrow
 H^{s-2}\left(
 M,T^*M\otimes\operatorname{ad}\mathbf{P}
 \right).
\]
The Yang--Mills gradient loses two derivatives; more precisely,
\begin{equation}
\label{eq:perdida-yang-mills}
 \mathcal Y\colon
 \mathcal A^s(\mathbf{P})
 \longrightarrow
 H^{s-2}\left(
 M,T^*M\otimes\operatorname{ad}\mathbf{P}
 \right),
 \qquad
 \mathcal Y(\mathbf{A}):=d_{\mathbf{A}}^*\mathbf{F}_{\mathbf{A}},
\end{equation}
is smooth.

\begin{remark}[The noncompact case]
On a connected noncompact manifold without boundary and of bounded geometry,
the compactness of the resolvent used in the preceding decomposition is no longer
available.
Theorem~\ref{teo:eichhorn-heber-slice-gauge} assumes that \(G\) is compact,
that the reference connection \(\mathbf{A}_0\) has bounded geometry and finite
Yang--Mills energy, and that
\[
 k-1\geq r>\frac{m}{2}+2.
\]
It also requires that the connection \(\mathbf{A}\) belong to the Sobolev component
under consideration and that
\[
 \inf\sigma_{\mathrm e}\left(
 \Delta_{\mathbf{A}}^{(0)}
 \restriction_{(\ker\Delta_{\mathbf{A}}^{(0)})^\perp}
 \right)>0.
\]
Under these hypotheses, one obtains a Hodge decomposition and a local
slice for the action of the reduced subgroup
\(\mathcal G_\bullet^{r+1}(\mathbf{A}_0)\).

The range in question is that of
$d_{\mathbf A}:H^{r+1}(M,\operatorname{ad}\mathbf P)\to
H^r(M,T^*M\otimes\operatorname{ad}\mathbf P)$; the complement is the
kernel of
$d_{\mathbf A}^*:H^r(M,T^*M\otimes\operatorname{ad}\mathbf P)
\to H^{r-1}(M,\operatorname{ad}\mathbf P)$. Orthogonality is measured
in $L^2(M,T^*M\otimes\operatorname{ad}\mathbf P)$, whereas the
sum and its projections are continuous in
$H^r(M,T^*M\otimes\operatorname{ad}\mathbf P)$.
The reduced subgroup is the one defined using the components of the full
group in the preceding chapter.
Proposition~\ref{prop:coulomb-sobolev-brecha-uniforme} constructs the
inverse on the complement of the kernel and the continuous projections.
Corollary~\ref{cor:carta-coulomb-sobolev-brecha} specifies the neighborhoods
for the local Coulomb parametrization.
Proposition~\ref{prop:tubo-calibre-brecha-conexa} extends that chart
to a tube over the orbit and proves that the orbit is closed.
\end{remark}

\section{Incompressible fluids as geodesics}
\label{sec:aplicacion-fluidos}
\index{Euler equations}
\index{incompressible fluid}
\index{diffeomorphism group!volume-preserving}
\index{right-invariant metric}

The third application uses a diffeomorphism group as its configuration
space. Incompressibility requires these
diffeomorphisms to preserve volume, and the kinetic energy defines a
weak right-invariant metric on this group. The corresponding geodesic
equation is the Euler equation.

Let \((M,\mathbf{g})\) be a compact, connected, oriented Riemannian
manifold of dimension $m$, with or without boundary, and let
\[
 s>\frac{m}{2}+1.
\]
Denote by \(\operatorname{Diff}^s_\mu(M)\) the identity component
of the subgroup of $\operatorname{Diff}^s(M)$ preserving
\(d\lambda_{\mathbf{g}}\) and \(\partial M\). When the boundary is empty,
conditions referring to it are omitted. In this section,
$\Delta_{\mathbf g}=\operatorname{div}_{\mathbf g}\nabla$ has the sign
of the Laplacian on functions from the chapter on Riemannian geometry; the
Hodge Laplacian $d\delta+\delta d$ used at the end is nonnegative.
Define
\begin{equation}
\label{eq:algebra-fluidos}
 \mathfrak X^s_\mu(M)
 :=
 \left\{
 \mathbf{u}\in H^s(M,TM)
 \middle|
 \operatorname{div}_{\mathbf{g}}\mathbf{u}=0,\quad
 (\mathbf{u}\cdot\boldsymbol{\nu})\restriction_{\partial M}=0
 \right\}.
\end{equation}
The same notation $\mathfrak X^t_\mu(M)$ will be used for any
$t>1/2$, interpreting the divergence distributionally and the normal
condition by means of the trace. In particular, it is defined for $t=s-1$.

\begin{lemma}[Global estimate for the Neumann problem]
\label{lem:neumann-regularidad-compacta-aplicaciones}
Let $(M,\mathbf g)$ be a smooth compact connected Riemannian manifold,
with or without boundary. For $a\geq0$, let $f\in H^a(M)$ and, when
$\partial M\neq\varnothing$, $h\in H^{a+1/2}(\partial M)$, with
\[
 \int_M f\,d\lambda_{\mathbf g}
 =\int_{\partial M}h\,d\sigma_{\mathbf g}.
\]
There exists a unique mean-zero solution $q\in H^{a+2}(M)$ of
$\Delta_{\mathbf g}q=f$, $\partial_{\boldsymbol\nu}q=h$, and
\begin{equation}
\label{eq:neumann-estimacion-datos-aplicaciones}
 \|q\|_{H^{a+2}(M)}
 \leq C_a\bigl(\|f\|_{H^a(M)}
                +\|h\|_{H^{a+1/2}(\partial M)}\bigr).
\end{equation}
When the boundary is empty, $h$ and its norm are omitted.
\end{lemma}

\begin{proof}
Begin with $h=0$ and mean-zero $f\in L^2(M)$.
Corollary~\ref{cor:poincare-promedio} makes the form
$\displaystyle \int_M\langle\nabla v,\nabla\phi\rangle_{\mathbf g}\,d\lambda_{\mathbf g}$
coercive on the mean-zero subspace of $H^1(M)$.
Theorem~\ref{teo:lax-milgram}, applied to the functional
$\phi\mapsto-\displaystyle \int_Mf\phi\,d\lambda_{\mathbf g}$, gives a weak solution
$v\in H^1(M)$ and
\begin{equation}
\label{eq:neumann-cota-debil-homogenea}
 \|v\|_{H^1(M)}\leq C\|f\|_{L^2(M)}.
\end{equation}
Since $f$ has mean zero, the weak identity holds for every
$\phi\in H^1(M)$, without any restriction on its mean or trace.

We now establish regularity up to the boundary. In an interior chart or
a chart straightening the boundary, write
$\omega=\sqrt{\det(\mathbf{g})}$ and $A^{ij}=\omega g^{ij}$. In a coordinate region
$Q$, which is a ball or a half-ball, the weak identity is
\begin{equation}
\label{eq:neumann-debil-coordenadas}
 \sum_{i,j=1}^m\int_Q A^{ij}\partial_jv\,\partial_i\phi\,dx
 =-\int_Q\omega f\phi\,dx.
\end{equation}
Test functions need not vanish on the flat face of the half-ball;
their supports remain away from the other faces. In a smaller chart,
the matrix $A$ is uniformly positive: $A^{ij}\xi_i\xi_j\geq
\lambda|\xi|^2$ for some $\lambda>0$.

Let $\chi$ be a cutoff supported away from those other faces, with
$\chi=1$ on a smaller region. For a coordinate direction $\alpha$,
tangential if there is a boundary, set
\[
 D_h^{\alpha} z(x)=\frac{z(x+h e_\alpha)-z(x)}h,
 \qquad
 \phi=-D_{-h}^{\alpha}(\chi^2D_h^{\alpha} v).
\]
For small $|h|$, this test function belongs to $H^1(Q)$ and is admissible.
For a function $w$ whose support and the translates in use stay in $Q$,
except possibly at its flat face, a change of variables gives
\[
 \int_Q z\,D_{-h}^{\alpha}w\,dx
 =-\int_Q(D_h^{\alpha}z)w\,dx.
\]
Apply this with $w=\chi^2D_h^{\alpha}v$; the tangential direction does
not cross the flat face. Moreover,
\[
 D_h^{\alpha}(A\nabla v)
 =A(\,\cdot+h e_\alpha)\nabla D_h^{\alpha} v
   +(D_h^{\alpha} A)\nabla v.
\]
The principal term in \eqref{eq:neumann-debil-coordenadas} is therefore
$\displaystyle \int_Q\chi^2 A^{ij}(x+h e_\alpha)
\partial_jD_h^{\alpha} v\,\partial_iD_h^{\alpha} v\,dx$.
The remaining terms contain $\nabla\chi$, $D_h^{\alpha} A$, or $f$.
We use
\[
 \|D_h^{\alpha} v\|_{L^2(\operatorname{supp}\chi)}
 \leq\|\partial_\alpha v\|_{L^2(Q)},\qquad
 \|D_h^{\alpha} A\|_{L^\infty(\operatorname{supp}\chi)}
 \leq\|\partial_\alpha A\|_{L^\infty(Q)}.
\]
For the term involving $f$, retain the difference on the test function:
\[
 \left|\int_Q\omega fD_{-h}^{\alpha}(\chi^2D_h^{\alpha} v)\,dx\right|
 \leq C\|f\|_{L^2(Q)}
       \bigl(\|\chi\nabla D_h^{\alpha} v\|_{L^2(Q,\mathbb R^m)}
                         +\|v\|_{H^1(Q)}\bigr).
\]
These inequalities follow from the integral formula for translations
and Cauchy--Schwarz; for functions in $H^1(Q)$, they follow by the
smooth approximation in Lemma~\ref{lem:calculo-sobolev-compacto-aplicaciones}(c), using a cutoff in a larger chart.
Applying Young's inequality to each product and absorbing half the principal
term yields
\begin{equation}
\label{eq:neumann-diferencias-tangenciales}
 \|\chi\nabla D_h^{\alpha} v\|_{L^2(Q,\mathbb R^m)}^2
 \leq C\bigl(\|f\|_{L^2(Q)}^2+\|v\|_{H^1(Q)}^2\bigr),
\end{equation}
with $C$ independent of $h$. When tested against compactly supported
smooth functions, a weakly convergent subsequence of
these differences has limit $\partial_\alpha\nabla v$. This gives in
$L^2(\Omega)$ all second-order derivatives except, in a boundary
chart, the double normal derivative.

The interior equation determines the remaining derivative:
\begin{equation}
\label{eq:neumann-recuperacion-normal}
 A^{mm}\partial_m^2v
 =\omega f-\sum_{\substack{1\leq i,j\leq m\\(i,j)\neq(m,m)}}A^{ij}\partial_i\partial_jv
          -\displaystyle\sum_{i,j=1}^m(\partial_iA^{ij})\partial_jv.
\end{equation}
All terms on the right belong to $L^2(\Omega)$, and $A^{mm}\geq\lambda$.
This distributional equality proves that $\partial_m^2v\in L^2(\Omega)$ on the
smaller region. In an interior chart, all difference directions are
allowed from the outset. A finite cover and the chosen cutoffs
give $v\in H^2(M)$ and
$\|v\|_{H^2(M)}\leq C(\|f\|_{L^2(M)}+\|v\|_{H^1(M)})$.

For higher orders, we proceed by induction. Suppose that
$f\in H^j(M)$ and that $v\in H^{j+1}(M)$ has already been obtained, with $j\geq1$.
If $\beta$ is a tangential multi-index of length $j$, differentiate
\eqref{eq:neumann-debil-coordenadas} tangentially in the distributional sense.
The field $w=\partial^\beta v\in H^1(Q)$ satisfies
\[
 \int_Q A\nabla w\cdot\nabla\phi\,dx
 =-\int_Q\partial^\beta(\omega f)\phi\,dx
   +\int_Q\mathbf G_\beta\cdot\nabla\phi\,dx,
\]
where
\[
 \mathbf G_\beta
 =-\displaystyle\sum_{0<\gamma\leq\beta}\binom\beta\gamma
       (\partial^\gamma A)\nabla\partial^{\beta-\gamma}v.
\]
Here $\partial^\beta(\omega f)\in L^2(Q)$ and
$\mathbf G_\beta\in H^1(Q,\mathbb R^m)$: each derivative of $v$ in
$\mathbf G_\beta$ has order at most $j$.
The same difference quotient proof applies to $w$. The additional term
becomes
$\displaystyle \int_QD_h^\alpha\mathbf G_\beta\cdot
\nabla(\chi^2D_h^\alpha w)\,dx$ and is estimated by
\[
 C\|\mathbf G_\beta\|_{H^1(Q,\mathbb R^m)}
 \bigl(\|\chi\nabla D_h^\alpha w\|_{L^2(Q,\mathbb R^m)}
                               +\|w\|_{H^1(Q)}\bigr).
\]
This gives all derivatives of order $j+2$ with at most one normal
index. For the others, differentiate
\eqref{eq:neumann-recuperacion-normal} $j$ times. A derivative with
$b\geq2$ normal indices is expressed through derivatives of the
same order with at most $b-1$ normal indices, derivatives of $v$ of
order at most $j+1$, and derivatives of $f$ of order at most $j$.
Induction on $b$ determines them all. Combining the finitely many
charts, induction on $j$ gives
\[
 \|v\|_{H^{j+2}(M)}
 \leq C_j\bigl(\|f\|_{H^j(M)}+\|v\|_{H^1(M)}\bigr),
 \qquad j\in\mathbb N_0.
\]
The preceding tangential differentiations are also justified by
successive differences; their distributional limits are the derivatives
written above, and the bounds use only orders already obtained in the induction.

For arbitrary $f\in L^2(M)$, let $Tf$ be the mean-zero weak
solution with datum $f-\overline f$, where
$\overline f=\lambda_{\mathbf g}(M)^{-1}\displaystyle \int_Mf\,d\lambda_{\mathbf g}$.
We have proved that the same operator $T$ is continuous from $H^j(M)$ to
$H^{j+2}(M)$ for each integer $j\geq0$. The interpolation in
Lemma~\ref{lem:calculo-sobolev-compacto-aplicaciones}, applied between
consecutive integers, gives
$T:H^a(M)\to H^{a+2}(M)$ for every $a\geq0$.
Its weak identity remains the one defining it in $L^2(M)$.

We now treat $h\neq0$.
Theorem~\ref{teo:traza-frontera-geometria-acotada}, with $p=2$ and a first-order
jet, provides $b\in H^{a+2}(M)$ with
\[
 b|_{\partial M}=0,\qquad \partial_{\boldsymbol\nu}b=h,
 \qquad
 \|b\|_{H^{a+2}(M)}\leq C_a\|h\|_{H^{a+1/2}(\partial M)}.
\]
The theorem is formulated using the inward collar coordinate: to
obtain the outward normal, prescribe the jet $(0,-h)$.
The divergence identity for $b\in H^2(M)$ follows by approximation
in $H^2(M)$ by smooth functions; the Laplacian converges in $L^2(M)$
and the normal derivative in $H^{1/2}(\partial M)$. Thus
$\displaystyle \int_M(f-\Delta_{\mathbf g}b)\,d\lambda_{\mathbf g}=0$.
The function
\[
 q=b+T(f-\Delta_{\mathbf g}b)-\overline b
\]
has the required mean, equation, and normal condition. The bounds
for $T$, $b$, and $\Delta_{\mathbf g}:H^{a+2}(M)\to H^a(M)$
prove \eqref{eq:neumann-estimacion-datos-aplicaciones}.
If two solutions have the same data, their difference $z$ satisfies
$\displaystyle \int_M|\nabla z|_{\mathbf g}^2\,d\lambda_{\mathbf g}=0$ upon using $z$
as a weak test function. Connectedness and the zero mean imply $z=0$.
\end{proof}

\begin{proposition}[Neumann problem and Hodge--Leray projection]
\label{prop:neumann-hodge-leray}
\index{Neumann problem}
\index{Hodge--Leray projection}
Let $(M,\mathbf g)$ be a smooth compact connected oriented Riemannian
manifold of dimension $m$, with or without boundary. For $r\geq0$ there exists a continuous linear
operator
\[
 Q:H^r(M,TM)\longrightarrow
 \left\{q\in H^{r+1}(M)\ \middle|\ \int_Mq\,d\lambda_{\mathbf g}=0\right\}
\]
characterized by
\begin{equation}
\label{eq:neumann-hodge-leray-debil}
 \int_M\langle\nabla Q\mathbf X,\nabla\phi\rangle_{\mathbf g}
                  \,d\lambda_{\mathbf g}
 =\int_M\langle\mathbf X,\nabla\phi\rangle_{\mathbf g}
                  \,d\lambda_{\mathbf g},\qquad \phi\in H^1(M).
\end{equation}
Set
\[
 \mathfrak X^0_\mu(M)
 :=\left\{\mathbf Y\in L^2(M,TM)\ \middle|\
 \int_M\langle\mathbf Y,\nabla\phi\rangle_{\mathbf g}
                 \,d\lambda_{\mathbf g}=0\quad\forall\phi\in H^1(M)\right\}.
\]
The operator $\mathbb P\mathbf X:=\mathbf X-\nabla Q\mathbf X$ is the
orthogonal projection of $L^2(M,TM)$ onto $\mathfrak X^0_\mu(M)$, and
\begin{equation}
\label{eq:cotas-Q-P-hodge-leray}
 \|Q\mathbf X\|_{H^{r+1}(M)}
 +\|\mathbb P\mathbf X\|_{H^r(M,TM)}
 \leq C_r\|\mathbf X\|_{H^r(M,TM)}.
\end{equation}
If $r>1/2$, its range in $H^r(M,TM)$ is $\mathfrak X^r_\mu(M)$, with
the definition of divergence and trace given above. In particular, for
$s>m/2+1$ and $\mathbf X\in H^{s-1}(M,TM)$, $q=Q\mathbf X\in H^s(M)$
solves
\begin{equation}
\label{eq:neumann-hodge-leray}
 \Delta_{\mathbf g}q=\operatorname{div}_{\mathbf g}\mathbf X,
 \qquad
 \partial_{\boldsymbol\nu}q
 =\langle\operatorname{Tr}\mathbf X,\boldsymbol\nu\rangle_{\mathbf g}.
\end{equation}
The first equality is distributional, and the second is an equality of
traces. We have the direct sum, orthogonal in $L^2(M,TM)$,
\[
 H^{s-1}(M,TM)=\mathfrak X^{s-1}_\mu(M)
           \oplus\nabla\bigl(H^s(M)/\mathbb R\bigr).
\]
Moreover, $Q$ and $\mathbb P$ map smooth data to data smooth up to the boundary.
\end{proposition}

\begin{proof}
For $\mathbf X\in L^2(M,TM)$, Poincaré and Lax--Milgram, used as
in the preceding lemma, apply to the functional
$\phi\mapsto\displaystyle \int_M\langle\mathbf X,\nabla\phi\rangle_{\mathbf g}
\,d\lambda_{\mathbf g}$. Its norm does not exceed
$\|\mathbf X\|_{L^2(M,TM)}$. This yields a unique mean-zero $Q\mathbf X$
with
$\|Q\mathbf X\|_{H^1(M)}\leq C\|\mathbf X\|_{L^2(M,TM)}$.
The weak identity directly proves that
$\mathbb P\mathbf X\in\mathfrak X^0_\mu(M)$ and that
$\mathbf X-\mathbb P\mathbf X$ is orthogonal to this subspace.
If $\mathbf X\in\mathfrak X^0_\mu(M)$, uniqueness gives $Q\mathbf X=0$;
thus $\mathbb P^2=\mathbb P$ and $\mathbb P$ is the asserted orthogonal
projection.

If $\mathbf X\in H^k(M,TM)$, with integer $k\geq1$, the data
$f=\operatorname{div}_{\mathbf g}\mathbf X$ and
$h=\langle\operatorname{Tr}\mathbf X,\boldsymbol\nu\rangle_{\mathbf g}$
belong to $H^{k-1}(M)$ and
$H^{k-1/2}(\partial M)$, respectively. To justify their compatibility, take
$\mathbf X_j\in\Gamma(TM)$ with
$\mathbf X_j\to\mathbf X$ in $H^1(M,TM)$, by
Lemma~\ref{lem:calculo-sobolev-compacto-aplicaciones}.
Theorem~\ref{teo:divergencia} holds for each $\mathbf X_j$.
Differentiation and the trace are continuous from $H^1(M,TM)$ to
$L^2(M)$ and $H^{1/2}(\partial M,TM|_{\partial M})$, respectively.
Moreover,
\[
 \left|\int_M\operatorname{div}_{\mathbf g}(\mathbf X_j-\mathbf X)
        \,d\lambda_{\mathbf g}\right|
 \leq\lambda_{\mathbf g}(M)^{1/2}
       \|\operatorname{div}_{\mathbf g}(\mathbf X_j-\mathbf X)\|_{L^2(M)},
\]
and the boundary integral satisfies the corresponding bound by
Cauchy--Schwarz on $\partial M$. Both errors tend to zero.
This proves the divergence formula for every vector field in $H^1(M,TM)$
and, in particular, compatibility of the specified data. The preceding lemma
with $a=k-1$ provides $q\in H^{k+1}(M)$ with
$\|q\|_{H^{k+1}(M)}\leq C_k\|\mathbf X\|_{H^k(M,TM)}$.
Integration by parts gives
\eqref{eq:neumann-hodge-leray-debil}, so $q=Q\mathbf X$ by
uniqueness. Interpolate this same operator between $k=0,1,2,\ldots$;
Lemma~\ref{lem:calculo-sobolev-compacto-aplicaciones} gives the bounds in
$Q$ for every $r\geq0$. Continuity of the gradient proves those for
$\mathbb P$.

Let us make the normal condition precise. If $r>1/2$ and
$\mathbf Y\in H^r(M,TM)$ has zero distributional divergence, then
\begin{equation}
\label{eq:green-divergencia-cero-traza}
 \int_M\langle\mathbf Y,\nabla\phi\rangle_{\mathbf g}\,d\lambda_{\mathbf g}
 =\int_{\partial M}\phi\langle\operatorname{Tr}\mathbf Y,
                  \boldsymbol\nu\rangle_{\mathbf g}\,d\sigma_{\mathbf g},
 \qquad\phi\in C^\infty(M).
\end{equation}
For $r\geq1$, this is Green's formula obtained by smooth approximation.
For $1/2<r<1$, the construction of the trace on the half-space implies
that restrictions to the collar hypersurfaces at distance $t$
converge in $L^2(\partial M,TM|_{\partial M})$ to the trace as
$t\to0^{+}$, identifying the fibers by means of the collar.
In the Euclidean proof, this convergence follows by Fourier analysis in the
tangential variables and Cauchy--Schwarz in the normal variable; the integral
$\displaystyle \int_{\mathbb R}(1+\xi_m^2)^{-r}\,d\xi_m$ is finite.
Take a cutoff $\zeta_\varepsilon$ that is zero near the boundary,
one outside the collar of width $\varepsilon$, and in the collar equals
$\zeta(t/\varepsilon)$, with $\zeta(0)=0$, $\zeta(1)=1$.
Applying the zero-divergence condition to $\phi\zeta_\varepsilon$ gives
\[
 0=\int_M\zeta_\varepsilon\langle\mathbf Y,\nabla\phi\rangle_{\mathbf g}
       \,d\lambda_{\mathbf g}
   +\int_M\phi\langle\mathbf Y,\nabla\zeta_\varepsilon\rangle_{\mathbf g}
       \,d\lambda_{\mathbf g}.
\]
The first term converges to the interior integral in
\eqref{eq:green-divergencia-cero-traza}. Convergence of the
collar restrictions and $\displaystyle \int_0^1\zeta'(t)\,dt=1$ make the second
term converge to the negative of the boundary integral, since the collar
coordinate increases inward. This proves the remaining identity.

If $\mathbf Y\in\mathfrak X^0_\mu(M)\cap H^r(M,TM)$, test functions
supported in the interior give $\operatorname{div}_{\mathbf g}\mathbf Y=0$.
The preceding identity and the extension of every smooth boundary function
give $\langle\operatorname{Tr}\mathbf Y,\boldsymbol\nu\rangle_{\mathbf g}=0$.
Conversely, these two conditions make the right-hand side of
\eqref{eq:green-divergencia-cero-traza} vanish; density of $C^\infty(M)$
in $H^1(M)$ extends the identity to all test functions in that space.
This identifies the range. Applied to
$\mathbf Y=\mathbf X-\nabla Q\mathbf X$, the same argument also gives
\eqref{eq:neumann-hodge-leray}. Finally, if $\mathbf X$ is smooth,
the bounds at each integer order give $Q\mathbf X\in H^{k+1}(M)$ for
every $k$. The Sobolev embeddings at each order prove that
$Q\mathbf X$ and $\mathbb P\mathbf X$ are smooth up to the boundary.
\end{proof}

\begin{corollary}[Approximation preserving incompressibility]
\label{cor:densidad-campos-solenoidales-compactos}
For $r\geq0$, smooth divergence-free vector fields tangent to
$\partial M$ are dense in
$\mathfrak X^0_\mu(M)\cap H^r(M,TM)$. In particular, they are dense in
$\mathfrak X^{s-1}_\mu(M)$ for $s>m/2+1$ and in
$\mathfrak X^0_\mu(M)$ with the $L^2(M,TM)$ norm.
\end{corollary}

\begin{proof}
Let $\mathbf Z$ be an element of the indicated space.
Lemma~\ref{lem:calculo-sobolev-compacto-aplicaciones} provides
$\mathbf V_j\in\Gamma(TM)$ with
$\|\mathbf V_j-\mathbf Z\|_{H^r(M,TM)}\to0$.
Using the mean-zero operator from the preceding proposition, define
\[
 q_j=Q\mathbf V_j,\qquad
 \mathbf Z_j=\mathbf V_j-\nabla q_j=\mathbb P\mathbf V_j.
\]
Each $q_j$ and each $\mathbf Z_j$ is smooth up to the boundary. The
Neumann equations give
$\operatorname{div}_{\mathbf g}\mathbf Z_j=0$ and
$\langle\mathbf Z_j,\boldsymbol\nu\rangle_{\mathbf g}=0$.
Since $\mathbb P\mathbf Z=\mathbf Z$, the explicit estimate is
\[
 \|\mathbf Z_j-\mathbf Z\|_{H^r(M,TM)}
 =\|\mathbb P(\mathbf V_j-\mathbf Z)\|_{H^r(M,TM)}
 \leq C_r\|\mathbf V_j-\mathbf Z\|_{H^r(M,TM)}\longrightarrow0.
\]
\end{proof}

\begin{theorem}[Ebin--Marsden: the volume-preserving group]
\label{teo:ebin-marsden-diffmu-compacta}
Let \((M,\mathbf{g})\) be a compact connected oriented Riemannian
manifold of dimension $m$, with or without boundary. If
\[
 s>\frac{m}{2}+1,
\]
then \(\operatorname{Diff}^s_\mu(M)\) is a closed \textit{split} Hilbert
submanifold of \(\operatorname{Diff}^s(M)\). For each
\(\eta\in\operatorname{Diff}^s_\mu(M)\), the map
\[
 \mathcal I_\eta\colon
 \mathfrak X^s_\mu(M)
 \longrightarrow
 T_\eta\operatorname{Diff}^s_\mu(M),
 \qquad
 \mathcal I_\eta(\mathbf{u}):=\mathbf{u}\circ\eta,
\]
is a topological linear isomorphism, with inverse
\[
 \mathbf{U}\longmapsto \mathbf{U}\circ\eta^{-1}.
\]
In particular,
\[
 T_{\operatorname{id}}\operatorname{Diff}^s_\mu(M)
 =
 \mathfrak X^s_\mu(M).
\]
\end{theorem}

\begin{proof}
Write $\mu=d\lambda_{\mathbf g}$ and define the Jacobian density
$J_\eta$ by $\eta^*\mu=J_\eta\mu$. Change of variables gives
$\displaystyle\int_MJ_\eta\,d\lambda_{\mathbf g}=\lambda_{\mathbf g}(M)$.
Thus the map
\[
 \mathcal J:\operatorname{Diff}^s(M)\longrightarrow Y^{s-1},
 \qquad \mathcal J(\eta)=J_\eta-1,
 \qquad Y^{s-1}:=
 \left\{f\in H^{s-1}(M)\mid\int_M f\,d\lambda_{\mathbf g}=0\right\},
\]
is well defined. In source and target charts,
\[
 J_\eta(x)=
 \frac{\sqrt{\det(\mathbf g(\eta(x)))}}{\sqrt{\det(\mathbf g(x))}}
                  |\det D\eta(x)|.
\]
In a $C^1$ neighborhood of any diffeomorphism, the determinant has
constant sign on each chart component and is bounded away from zero.
Smooth superposition and the algebra property of $H^{s-1}(M)$, since
$s-1>m/2$, prove that $\mathcal J$ is smooth. In boundary charts the
expressions are restrictions of smooth local extensions of the metric
coefficients. The charts of $\operatorname{Diff}^s(M)$ already incorporate
preservation of the boundary.

At the identity, differentiation of the determinant and metric density gives
\[
 D\mathcal J_{\operatorname{id}}\mathbf X
 =\frac{1}{\sqrt{\det\mathbf g}}\sum_{i=1}^m
       \partial_i\bigl(\sqrt{\det\mathbf g}\,X^i\bigr)
 =\operatorname{div}_{\mathbf g}\mathbf X.
\]
The identity is first obtained for smooth fields tangent to the boundary
and extends continuously to $H^s(M,TM)$. It also follows pointwise from the
differential of the Jacobian formula, valid for $C^1$ fields. Its integral
vanishes by the divergence formula and the zero normal component.

For $f\in Y^{s-1}$, Lemma~\ref{lem:neumann-regularidad-compacta-aplicaciones},
with $a=s-1$, gives a unique mean-zero function $q_f\in H^{s+1}(M)$ satisfying
\[
 \Delta_{\mathbf g}q_f=f,\qquad
 \partial_{\boldsymbol\nu}q_f=0\text{ on }\partial M,\qquad
 \|q_f\|_{H^{s+1}(M)}\leq C_N\|f\|_{H^{s-1}(M)}.
\]
If the boundary is empty, omit the normal condition. Thus $R f:=\nabla q_f$
belongs to the tangent space of $\operatorname{Diff}^s(M)$ at the identity,
depends linearly and continuously on $f$, and satisfies
$D\mathcal J_{\operatorname{id}}R f=f$. In particular,
$R D\mathcal J_{\operatorname{id}}$ is a continuous projection and
\[
 T_{\operatorname{id}}\operatorname{Diff}^s(M)
 =\mathfrak X^s_\mu(M)\oplus R Y^{s-1}.
\]
The Banach implicit function theorem makes $\mathcal J^{-1}(0)$ a split
submanifold near the identity, with tangent space $\mathfrak X^s_\mu(M)$.

Now let $\eta\in\mathcal H:=\mathcal J^{-1}(0)$. Right translation
$R_\eta:\varphi\mapsto\varphi\circ\eta$ is a Hilbert diffeomorphism:
in charts it composes with a fixed diffeomorphism, and its inverse is
$R_{\eta^{-1}}$. Moreover,
$J_{\varphi\circ\eta}=(J_\varphi\circ\eta)J_\eta=J_\varphi\circ\eta$,
so it preserves the condition $\mathcal J=0$. It therefore transports
the preceding local submanifold to a neighborhood of $\eta$; its
differential is $\mathbf X\mapsto\mathbf X\circ\eta$. The composition
bounds for $\eta$ and $\eta^{-1}$ make this tangent identification a
linear topological isomorphism. Continuity of $\mathcal J$ makes $\mathcal H$ closed. Since $\mathcal H$
is a Hilbert manifold, it is locally path connected; its identity
component is open in $\mathcal H$ and, as every connected component,
is also closed. Thus $\operatorname{Diff}^s_\mu(M)$ is a closed split
submanifold of $\operatorname{Diff}^s(M)$, with the same tangent spaces as $\mathcal H$.
\end{proof}

The weak right-invariant kinetic metric is
\begin{equation}
\label{eq:metrica-cinetica-diffmu}
 \mathbf{G}_\eta(\mathbf{U},\mathbf{V})
 =
 \int_M
 \left\langle
 \mathbf{U}\circ\eta^{-1},
 \mathbf{V}\circ\eta^{-1}
 \right\rangle_{\mathbf{g}}
 \,d\lambda_{\mathbf{g}}.
\end{equation}
Since $\eta$ preserves volume, a change of variables transforms this
formula into
\[
 \mathbf G_\eta(\mathbf U,\mathbf V)
 =\int_M \mathbf g_{\eta(x)}(\mathbf U(x),\mathbf V(x))\,d\lambda_{\mathbf g}(x).
\]
Smooth dependence of this metric is verified in the latter expression
by superposition and products, without assuming smoothness of inversion on
$\operatorname{Diff}^s(M)$. For a curve $\eta(t)$, define its Eulerian
velocity by
\[
 \mathbf{u}(t):=\dot\eta(t)\circ\eta(t)^{-1}.
\]
The kinetic action is
\begin{equation}
\label{eq:accion-cinetica-fluidos}
 \mathcal A(\eta)
 :=
 \frac{1}{2}
 \int_{t_0}^{t_1}
 \int_M|\mathbf{u}(t,x)|_{\mathbf{g}}^2\,d\lambda_{\mathbf{g}}\,dt.
\end{equation}

\begin{proposition}[Euler as a geodesic equation]
\label{prop:euler-geodesica}
Let $(M,\mathbf g)$ be a smooth compact connected oriented Riemannian
manifold of dimension $m$, with or without boundary, and let $s>m/2+1$. A curve
$\eta\in C^2([t_0,t_1],\operatorname{Diff}^s_\mu(M))$, with fixed
endpoints, is critical for \eqref{eq:accion-cinetica-fluidos} if and only if its
Eulerian velocity satisfies
\begin{equation}
\label{eq:euler-incompresible}
 \partial_t\mathbf u+\nabla_{\mathbf u}\mathbf u=-\nabla p,
 \qquad \operatorname{div}_{\mathbf g}\mathbf u=0,
 \qquad
 \langle\operatorname{Tr}\mathbf u,\boldsymbol\nu\rangle_{\mathbf g}=0
 \quad\text{on }\partial M.
\end{equation}
Here
\[
 \mathbf u\in C([t_0,t_1],H^s(M,TM))
       \cap C^1([t_0,t_1],H^{s-1}(M,TM)),
 \qquad p\in C([t_0,t_1],H^s(M)),
\]
and the first equation is interpreted in $H^{s-1}(M,TM)$.
The pressure is determined at each time up to a spatial constant;
the normalization
$\displaystyle \int_Mp(t,x)\,d\lambda_{\mathbf g}(x)=0$ makes it unique.
The boundary condition is omitted when $\partial M=\varnothing$.
\end{proposition}

\begin{proof}
Write $I=[t_0,t_1]$. We use variations of class $C^2$
in $(\varepsilon,t)$ with values in $\operatorname{Diff}^s_\mu(M)$
and fixed endpoints. Continuous composition from
Lemma~\ref{lem:composicion-niveles-inferiores-curvas} and continuous
inversion from Theorem~\ref{teo:ebin-marsden-diff-sobolev} give
$\mathbf u=\dot\eta\circ\eta^{-1}\in C(I,H^s(M,TM))$.
The embedding $H^s(M,TM)\hookrightarrow C^1(M,TM)$ allows us to differentiate
the identity $\dot\eta=\mathbf u\circ\eta$ pointwise. This gives
\begin{equation}
\label{eq:aceleracion-euleriana-lagrangiana}
 \partial_t\mathbf u
 =(\nabla_t\dot\eta)\circ\eta^{-1}-\nabla_{\mathbf u}\mathbf u.
\end{equation}
The covariant acceleration $\nabla_t\dot\eta$ is continuous in the space
of sections $H^s(M,\eta(t)^*TM)$ along $\eta$: in coordinates, its
components are
$\ddot\eta^a+\displaystyle\sum_{b,c=1}^{m}\Gamma^a_{bc}(\eta)\dot\eta^b\dot\eta^c$.
Superposition and the product estimate in
Lemma~\ref{lem:calculo-sobolev-compacto-aplicaciones} control each
term. By \eqref{eq:estimacion-transporte-sobolev-aplicaciones}, the
right-hand side of \eqref{eq:aceleracion-euleriana-lagrangiana} is
continuous in $H^{s-1}(M,TM)$. The last part of
Lemma~\ref{lem:composicion-niveles-inferiores-curvas}, or its proof
using the Bochner integral of this right-hand side, gives
$\mathbf u\in C^1(I,H^{s-1}(M,TM))$.
Thus this regularity follows from the given curve, before imposing the
Euler equation.

Let $\eta_\varepsilon$ be an admissible variation, and set
\[
 \mathbf W(t)=\left.\partial_\varepsilon\eta_\varepsilon(t)\right|_0,
 \qquad \mathbf w=\mathbf W\circ\eta^{-1}.
\]
The tangent-space identification in
Theorem~\ref{teo:ebin-marsden-diffmu-compacta} implies
$\mathbf w(t)\in\mathfrak X^s_\mu(M)$ and
$\mathbf w(t_0)=\mathbf w(t_1)=0$.
The same argument applied to $\mathbf W$ shows that
\[
 \mathbf w\in C(I,H^s(M,TM))\cap C^1(I,H^{s-1}(M,TM)),
 \quad
 \partial_t\mathbf w
 =(\nabla_t\mathbf W)\circ\eta^{-1}-\nabla_{\mathbf u}\mathbf w.
\]

Conversely, let $\mathbf w$ be smooth in $I\times M$, divergence-free,
tangent to the boundary, and zero at the temporal endpoints. For each
$t$, take the flow $\phi_\varepsilon^t$ of the spatial vector field
$\mathbf w(t,\cdot)$. Theorem~\ref{teo: fundamental de flujos}
and the dependence of ODEs on parameters give a
common interval in $\varepsilon$ and smooth dependence on $(\varepsilon,t,x)$;
a common choice is possible because $I\times M$ is compact.
Tangency makes the boundary invariant by uniqueness of the flow, and
Lie's formula for the volume form gives
\[
 \frac{d}{d\varepsilon}(\phi_\varepsilon^t)^*d\lambda_{\mathbf g}
 =(\phi_\varepsilon^t)^*
       ((\operatorname{div}_{\mathbf g}\mathbf w(t))d\lambda_{\mathbf g})=0.
\]
Moreover, $\mathbf w(t_i)=0$ implies
$\phi_\varepsilon^{t_i}=\operatorname{id}$ for $i=0,1$.
By smooth superposition,
\[
 \eta_\varepsilon(t)=\phi_\varepsilon^t\circ\eta(t)
\]
is an admissible variation of class $C^2$ with Eulerian vector field
$\mathbf w$. The path in $\varepsilon$ starting at $\eta(t)$ ensures
that it remains in the identity component.

We now differentiate the velocity of a general variation. The version of the
preceding calculation with parameter $\varepsilon$ proves the existence of
$\delta\mathbf u=\partial_\varepsilon\mathbf u_\varepsilon|_0$ in
$C(I,H^{s-1}(M,TM))$. Differentiate
$\dot\eta_\varepsilon=\mathbf u_\varepsilon\circ\eta_\varepsilon$ and
$\mathbf W=\mathbf w\circ\eta$ covariantly. This gives
\[
 \left.\nabla_\varepsilon\dot\eta_\varepsilon\right|_0
 =(\delta\mathbf u+\nabla_{\mathbf w}\mathbf u)\circ\eta,
 \qquad
 \nabla_t\mathbf W
 =(\partial_t\mathbf w+\nabla_{\mathbf u}\mathbf w)\circ\eta.
\]
The Levi--Civita connection is torsion-free, and the mixed derivatives of
$\eta_\varepsilon$ commute. The two left-hand sides are equal,
and we obtain Lin's identity
\begin{equation}
\label{eq:identidad-lin-sobolev}
 \delta\mathbf u=\partial_t\mathbf w+[\mathbf u,\mathbf w]
 =\partial_t\mathbf w+\nabla_{\mathbf u}\mathbf w
                         -\nabla_{\mathbf w}\mathbf u
 \quad\text{in }C(I,H^{s-1}(M,TM)).
\end{equation}
All contractions belong to this space by
\eqref{eq:estimacion-transporte-sobolev-aplicaciones}.

The functional $\mathbf v\mapsto\frac12\|\mathbf v\|_{L^2(M,TM)}^2$
is differentiable, and $H^{s-1}(M,TM)\hookrightarrow L^2(M,TM)$.
The preceding differentiation, uniform on the compact interval $I$,
allows differentiation under the time integral:
\[
 \delta\mathcal A
 =\int_I\int_M\langle\mathbf u,\delta\mathbf u\rangle_{\mathbf g}
                 \,d\lambda_{\mathbf g}\,dt.
\]
Since $\mathbf u$ and $\mathbf w$ are curves of class $C^1$ in
$L^2(M,TM)$, the product rule in that space gives
\[
 \int_I\int_M\langle\mathbf u,\partial_t\mathbf w\rangle_{\mathbf g}
                     \,d\lambda_{\mathbf g}\,dt
 =-\int_I\int_M\langle\partial_t\mathbf u,\mathbf w\rangle_{\mathbf g}
                     \,d\lambda_{\mathbf g}\,dt.
\]
The endpoint term is zero because $\mathbf w(t_i)=0$.
For spatial integration, we directly use
$\mathbf u(t),\mathbf w(t)\in C^1(M,TM)$, metric compatibility, and
the divergence formula in $H^1(M,TM)$ justified in
Proposition~\ref{prop:neumann-hodge-leray}. The vector fields to which it is
applied even belong to $H^s(M,TM)$ by the product estimate.
Indeed,
\begin{align*}
 \operatorname{div}_{\mathbf g}
       (\langle\mathbf u,\mathbf w\rangle_{\mathbf g}\mathbf u)
 &=\langle\nabla_{\mathbf u}\mathbf u,\mathbf w\rangle_{\mathbf g}
   +\langle\mathbf u,\nabla_{\mathbf u}\mathbf w\rangle_{\mathbf g},\\
 \operatorname{div}_{\mathbf g}
       (\tfrac12|\mathbf u|_{\mathbf g}^2\mathbf w)
 &=\langle\mathbf u,\nabla_{\mathbf w}\mathbf u\rangle_{\mathbf g}.
\end{align*}
We have used $\operatorname{div}_{\mathbf g}\mathbf u=
\operatorname{div}_{\mathbf g}\mathbf w=0$.
The normal fluxes of both displayed vector fields are zero. Consequently,
\begin{align*}
 \int_M\langle\mathbf u,\nabla_{\mathbf u}\mathbf w\rangle_{\mathbf g}
                    \,d\lambda_{\mathbf g}
 &=-\int_M\langle\nabla_{\mathbf u}\mathbf u,\mathbf w\rangle_{\mathbf g}
                    \,d\lambda_{\mathbf g},\\
 -\int_M\langle\mathbf u,\nabla_{\mathbf w}\mathbf u\rangle_{\mathbf g}
                    \,d\lambda_{\mathbf g}
 &=-\tfrac12\int_M\mathbf w(|\mathbf u|_{\mathbf g}^2)
                    \,d\lambda_{\mathbf g}=0.
\end{align*}
Thus
\begin{equation}
\label{eq:primera-variacion-euler-sobolev}
 \delta\mathcal A
 =-\int_I\int_M\langle\mathbf F(t),\mathbf w(t)\rangle_{\mathbf g}
                         \,d\lambda_{\mathbf g}\,dt,
 \qquad
 \mathbf F=\partial_t\mathbf u+\nabla_{\mathbf u}\mathbf u
       \in C(I,H^{s-1}(M,TM)).
\end{equation}

Suppose that $\eta$ is critical. For a smooth vector field
$\mathbf Z\in\mathfrak X^s_\mu(M)$ and
$\chi\in C_c^\infty((t_0,t_1))$, the vector field
$\mathbf w(t,x)=\chi(t)\mathbf Z(x)$ arises from the variation constructed
above. By \eqref{eq:primera-variacion-euler-sobolev},
\[
 \int_{t_0}^{t_1}\chi(t)
       \langle\mathbf F(t),\mathbf Z\rangle_{L^2(M,TM)}\,dt=0.
\]
The scalar function obtained by pairing is continuous. The fundamental lemma of the
calculus of variations makes it vanish on $(t_0,t_1)$, and continuity makes it vanish
at the endpoints as well. Fix $t\in I$. By
Corollary~\ref{cor:densidad-campos-solenoidales-compactos}, there exist smooth
vector fields $\mathbf Z_j\in\mathfrak X^s_\mu(M)$ with
\[
 \mathbf Z_j\longrightarrow\mathbb P\mathbf F(t)
 \quad\text{in }H^{s-1}(M,TM).
\]
By Cauchy--Schwarz and orthogonality of $\mathbb P$,
\[
 0=\lim_{j\to\infty}\langle\mathbf F(t),\mathbf Z_j\rangle_{L^2(M,TM)}
  =\langle\mathbf F(t),\mathbb P\mathbf F(t)\rangle_{L^2(M,TM)}
  =\|\mathbb P\mathbf F(t)\|_{L^2(M,TM)}^2.
\]
Thus $\mathbf F(t)=\nabla Q\mathbf F(t)$. Define
$p(t)=-Q\mathbf F(t)$. Proposition~\ref{prop:neumann-hodge-leray}
gives the equation, the zero mean, and, for $t,\tau\in I$,
\[
 \|p(t)-p(\tau)\|_{H^s(M)}
 \leq C_{s-1}\|\mathbf F(t)-\mathbf F(\tau)\|_{H^{s-1}(M,TM)}.
\]
This proves continuity of the pressure in time. If two pressures
produce the same equation, their difference has zero gradient; by
connectedness it is constant, and its zero mean makes it vanish.

Conversely, suppose that \eqref{eq:euler-incompresible}. For every
admissible variation, $\mathbf w(t)\in\mathfrak X^0_\mu(M)$ and
$p(t)\in H^s(M)\subseteq H^1(M)$. The definition of
$\mathfrak X^0_\mu(M)$ implies
$\displaystyle \int_M\langle\nabla p(t),\mathbf w(t)\rangle_{\mathbf g}
\,d\lambda_{\mathbf g}=0$. Substituting $\mathbf F=-\nabla p$ into
\eqref{eq:primera-variacion-euler-sobolev} gives $\delta\mathcal A=0$.
\end{proof}

The subspace $\mathfrak X^0_\mu(M)$ is closed in $L^2(M,TM)$,
since it is the intersection of the kernels of the continuous functionals
$\mathbf v\mapsto\langle\mathbf v,\nabla\phi\rangle_{L^2(M,TM)}$,
with $\phi\in H^1(M)$. Since $\mathbf u(t)$ belongs to this subspace,
so do its difference quotients. Their convergence
to $\partial_t\mathbf u(t)$ in $H^{s-1}(M,TM)$, and hence in
$L^2(M,TM)$, gives
$\partial_t\mathbf u(t)\in\mathfrak X^{s-1}_\mu(M)$.
Thus $\mathbb P\partial_t\mathbf u=\partial_t\mathbf u$ and
$Q\partial_t\mathbf u=0$. Applying $\mathbb P$ to
\eqref{eq:euler-incompresible} eliminates the pressure gradient and gives
the reduced equation
\begin{equation}
\label{eq:euler-proyectada}
 \partial_t\mathbf{u}=-\mathbb P(\nabla_{\mathbf{u}}\mathbf{u}).
\end{equation}
For \(\mathbf{X}=\nabla_{\mathbf{u}}\mathbf{u}\), let \(q\) be the function from
Proposition~\ref{prop:neumann-hodge-leray}. The identity
\[
 \nabla_{\mathbf{u}}\mathbf{u}
 =
 \mathbb P(\nabla_{\mathbf{u}}\mathbf{u})+\nabla q,
\],
substituted into \eqref{eq:euler-proyectada}, shows that \(p=-q\) and that
\[
 \Delta_{\mathbf{g}}p
 =
 -\operatorname{div}_{\mathbf{g}}(\nabla_{\mathbf{u}}\mathbf{u}),
 \qquad
 \partial_{\boldsymbol{\nu}}p
 =
 -\langle\nabla_{\mathbf{u}}\mathbf{u},\boldsymbol{\nu}\rangle_{\mathbf{g}}.
\]

In Eulerian coordinates, for the fixed metric \(\mathbf{g}\) and the boundary
conditions incorporated in \(\mathbb P_{\mathbf{g}}\), the right-hand side defines the
smooth quadratic map
\[
 \mathscr E_{\mathbf{g}}\colon
 \mathfrak X^s_\mu(M)
 \longrightarrow
 \mathfrak X^{s-1}_\mu(M),
 \qquad
 \mathscr E_{\mathbf{g}}(\mathbf{u})
 :=
 -\mathbb P_{\mathbf{g}}(\nabla_{\mathbf{u}}\mathbf{u}).
\]
Since
\[
 s-1>\frac{m}{2},
\]
the Sobolev product estimate gives
\[
 \|\nabla_{\mathbf{u}}\mathbf{u}\|_{H^{s-1}(M,TM)}
 \leq
 C_{M,\mathbf{g},s}\|\mathbf{u}\|_{H^s(M,TM)}^2.
\]
Proposition~\ref{prop:neumann-hodge-leray} proves continuity of
\(\mathbb P_{\mathbf{g}}\) on \(H^{s-1}(M,TM)\). Thus the Eulerian formulation
loses one derivative. Its linearization along a solution is
\begin{equation}
\label{eq:euler-linealizada}
 \partial_t\mathbf{v}+\nabla_{\mathbf{u}}\mathbf{v}
 +\nabla_{\mathbf{v}}\mathbf{u}=-\nabla q,
 \qquad
 \operatorname{div}_{\mathbf{g}}\mathbf{v}=0,
 \qquad
 (\mathbf{v}\cdot\boldsymbol{\nu})\restriction_{\partial M}=0.
\end{equation}

\begin{proposition}[Regularity of the pressure]
\label{prop:presion-ganancia-sobolev-compacta}
Under the hypotheses of Proposition~\ref{prop:euler-geodesica}, for
each $\mathbf u\in\mathfrak X^s_\mu(M)$ the mean-zero function
\[
 p(\mathbf u)=-Q(\nabla_{\mathbf u}\mathbf u)
\]
belongs to $H^{s+1}(M)$. If
$\operatorname{II}_{\mathrm{ext}}(\mathbf V,\mathbf W)
=\langle\nabla_{\mathbf V}\boldsymbol\nu,\mathbf W\rangle_{\mathbf g}$
for vectors tangent to the boundary, then
\begin{align}
 \Delta_{\mathbf g}p(\mathbf u)
 &=-\sum_{i,j=1}^m(\nabla_i u^j)(\nabla_j u^i)
                  -\operatorname{Ric}_{\mathbf g}(\mathbf u,\mathbf u),
 \label{eq:presion-dato-interior-regular}\\
 \partial_{\boldsymbol\nu}p(\mathbf u)
 &=\operatorname{II}_{\mathrm{ext}}
          (\operatorname{Tr}\mathbf u,\operatorname{Tr}\mathbf u).
 \label{eq:presion-dato-frontera-regular}
\end{align}
The interior datum belongs to $H^{s-1}(M)$ and the boundary datum to
$H^{s-1/2}(\partial M)$. There exist constants $C_s,C'_s$ such that
\begin{align}
 \|p(\mathbf u)\|_{H^{s+1}(M)}
 +\|\nabla p(\mathbf u)\|_{H^s(M,TM)}
 &\leq C_s\|\mathbf u\|_{H^s(M,TM)}^2,
 \label{eq:presion-cota-ganancia}\\
 \|p(\mathbf u)-p(\mathbf v)\|_{H^{s+1}(M)}
 &\leq C'_s\bigl(\|\mathbf u\|_{H^s(M,TM)}
                 +\|\mathbf v\|_{H^s(M,TM)}\bigr)
                   \|\mathbf u-\mathbf v\|_{H^s(M,TM)}.
 \label{eq:presion-continuidad-ganancia}
\end{align}
In particular, the normalized pressure of an Euler solution belongs
to $C(I,H^{s+1}(M))$.
\end{proposition}

\begin{proof}
Begin with a smooth divergence-free vector field tangent to the
boundary. In an orthonormal frame normal at the point of calculation,
write $\nabla_i u^j=\langle\nabla_{\mathbf e_i}\mathbf u,
\mathbf e_j\rangle_{\mathbf g}$. The Leibniz rule gives
\[
 \operatorname{div}_{\mathbf g}(\nabla_{\mathbf u}\mathbf u)
 =\sum_{i,j=1}^m(\nabla_i u^j)(\nabla_j u^i)
                 +\sum_{i,j=1}^m u^j\nabla_i\nabla_j u^i.
\]
The definition of curvature as a commutator implies
\[
 \sum_{i=1}^{m}\nabla_i\nabla_j u^i
 =\nabla_j\left(\sum_{i=1}^{m}\nabla_i u^i\right)
                   +\sum_{k=1}^{m}\operatorname{Ric}_{kj}u^k.
\]
Multiplying by $u^j$ and summing, the first part vanishes because
$\operatorname{div}_{\mathbf g}\mathbf u=0$. We obtain
\[
 \operatorname{div}_{\mathbf g}(\nabla_{\mathbf u}\mathbf u)
 =\sum_{i,j=1}^m(\nabla_i u^j)(\nabla_j u^i)
                      +\operatorname{Ric}_{\mathbf g}(\mathbf u,\mathbf u).
\]
On the boundary, differentiate
$\langle\mathbf u,\boldsymbol\nu\rangle_{\mathbf g}=0$ tangentially in the
direction $\mathbf u$ and use metric compatibility:
\[
 0=\langle\nabla_{\mathbf u}\mathbf u,\boldsymbol\nu\rangle_{\mathbf g}
                   +\langle\mathbf u,\nabla_{\mathbf u}\boldsymbol\nu\rangle_{\mathbf g}.
\]
This gives both pressure data in the smooth case.

For a vector field $\mathbf u\in\mathfrak X^s_\mu(M)$,
Corollary~\ref{cor:densidad-campos-solenoidales-compactos} gives smooth
vector fields $\mathbf u_j\in\mathfrak X^s_\mu(M)$ converging in
$H^s(M,TM)$. Let $f(\mathbf u)$ and $h(\mathbf u)$ be the right-hand
sides of \eqref{eq:presion-dato-interior-regular} and
\eqref{eq:presion-dato-frontera-regular}. The contraction defining
$f$ is independent of the frame, as the smooth identity shows.
Lemma~\ref{lem:calculo-sobolev-compacto-aplicaciones}, with
$s-1>m/2$, gives
\begin{align*}
 \|f(\mathbf u)\|_{H^{s-1}(M)}
 &\leq C_s\|\mathbf u\|_{H^s(M,TM)}^2,\\
 \|f(\mathbf u_j)-f(\mathbf u)\|_{H^{s-1}(M)}
 &\leq C_s\bigl(\|\mathbf u_j\|_{H^s(M,TM)}
                 +\|\mathbf u\|_{H^s(M,TM)}\bigr)
                  \|\mathbf u_j-\mathbf u\|_{H^s(M,TM)}.
\end{align*}
To verify the second bound, simply decompose each difference of
products as $(A_j-A)B_j+A(B_j-B)$ and use the first product
estimate. On the other hand,
$\nabla_{\mathbf u_j}\mathbf u_j\to\nabla_{\mathbf u}\mathbf u$
in $H^{s-1}(M,TM)$, so their divergences converge in
$H^{s-2}(M)$. The smooth identity and the stronger convergence of the
$f(\mathbf u_j)$ identify, in the distributional sense,
$-\operatorname{div}_{\mathbf g}(\nabla_{\mathbf u}\mathbf u)=f(\mathbf u)$.

The continuous trace from
Theorem~\ref{teo:traza-haz-frontera-geometria-acotada} gives
\[
 \operatorname{Tr}\mathbf u_j\to\operatorname{Tr}\mathbf u
 \quad\text{in }H^{s-1/2}(\partial M,TM|_{\partial M}).
\]
Since $s-1/2>(m-1)/2$, the same product lemma, now on
$\partial M$, yields
\begin{align*}
 \|h(\mathbf u)\|_{H^{s-1/2}(\partial M)}
 &\leq C_s\|\mathbf u\|_{H^s(M,TM)}^2,\\
 \|h(\mathbf u_j)-h(\mathbf u)\|_{H^{s-1/2}(\partial M)}
 &\leq C_s\bigl(\|\mathbf u_j\|_{H^s(M,TM)}
                 +\|\mathbf u\|_{H^s(M,TM)}\bigr)
                   \|\mathbf u_j-\mathbf u\|_{H^s(M,TM)}.
\end{align*}
The trace of $\nabla_{\mathbf u_j}\mathbf u_j$ converges in
$H^{s-3/2}(\partial M,TM|_{\partial M})$, since $s-1>1/2$.
Taking limits in the smooth boundary identity gives
$-\langle\operatorname{Tr}(\nabla_{\mathbf u}\mathbf u),
\boldsymbol\nu\rangle_{\mathbf g}=h(\mathbf u)$, with the higher regularity
just proved. If $m=1$, the boundary is finite and its spaces
of sections are finite-dimensional; the same conclusion is immediate.

For each $j$, the divergence formula gives
$\displaystyle \int_M f(\mathbf u_j)\,d\lambda_{\mathbf g}
=\int_{\partial M}h(\mathbf u_j)\,d\sigma_{\mathbf g}$.
The preceding convergences allow us to take limits in both integrals.
Thus the data $f(\mathbf u),h(\mathbf u)$ satisfy Neumann
compatibility.
Lemma~\ref{lem:neumann-regularidad-compacta-aplicaciones}, with $a=s-1$,
gives a mean-zero solution in $H^{s+1}(M)$ and the bound
\eqref{eq:presion-cota-ganancia}. By weak uniqueness, it agrees with
$-Q(\nabla_{\mathbf u}\mathbf u)$.
Applying the same Neumann estimate to differences of the data,
with the product bounds already written, proves
\eqref{eq:presion-continuidad-ganancia}. This last inequality,
applied to $\mathbf u(t)$ and $\mathbf u(\tau)$, proves continuity
in time in $H^{s+1}(M)$. Finally,
$\nabla_t\dot\eta=-(\nabla p)\circ\eta$ by
\eqref{eq:aceleracion-euleriana-lagrangiana}, and the Lagrangian
acceleration takes values at the $H^s(M,\eta^*TM)$ level of sections.
\end{proof}

To study dependence on the diffeomorphism, we use divergence, the exterior
derivative and, when a boundary is present, the normal component.
All Sobolev spaces of forms below are defined using the fixed metric
$\mathbf g$.

\begin{lemma}[Control of a vector field by its first-order data]
\label{lem:control-primer-orden-hodge-fluidos}
Let $(M,\mathbf g)$ be smooth, compact and connected, with or without boundary.
Let $\mathcal H$ denote the space of smooth vector fields $Z$ such that
\[
 d(Z^\flat)=0,\qquad \operatorname{div}_{\mathbf g}Z=0,
 \qquad \langle Z,\boldsymbol\nu\rangle_{\mathbf g}=0
                    \quad\text{on }\partial M.
\]
The last condition is omitted when the boundary is empty.
The space $\mathcal H$ is finite-dimensional. Let $\Pi_{\mathcal H}$
be its orthogonal projection in $L^2(M,TM)$. For every $s\geq1$ there is
$C_s>0$ such that
\begin{equation}
\label{eq:estimacion-primer-orden-hodge-fluidos}
 \begin{split}
 \|W\|_{H^s(M,TM)}\leq C_s\bigl(&
    \|d(W^\flat)\|_{H^{s-1}(M,\Lambda^2T^*M)}
   +\|\operatorname{div}_{\mathbf g}W\|_{H^{s-1}(M)}\\
   &+\|\langle\operatorname{Tr}W,\boldsymbol\nu\rangle_{\mathbf g}
                         \|_{H^{s-1/2}(\partial M)}
   +\|\Pi_{\mathcal H}W\|_{L^2(M,TM)}\bigr).
 \end{split}
\end{equation}
The boundary norm is omitted when $\partial M=\varnothing$.
\end{lemma}

\begin{proof}
The elliptic Hodge estimate with prescribed normal component is
\begin{equation}
\label{eq:hodge-primer-orden-con-resto-L2}
 \|\alpha\|_{H^s(M,T^*M)}
 \leq C_s^0\bigl(\|d\alpha\|_{H^{s-1}(M,\Lambda^2T^*M)}+
       \|\delta\alpha\|_{H^{s-1}(M)}+
       \|\alpha(\boldsymbol\nu)\|_{H^{s-1/2}(\partial M)}+
       \|\alpha\|_{L^2(M,T^*M)}\bigr).
\end{equation}
The interior norms are taken over $M$: $\alpha$ is a section
of $T^*M$, $d\alpha$ of $\Lambda^2T^*M$, and $\delta\alpha$ is scalar.
This is the regularity estimate for the Hodge problem used in
\cite[\S7, pp.~120--122]{EbinMarsden1970}; it holds for real orders
$s\geq1$ and has the corresponding weak regularity statement.

Let us specify the boundary conditions involved. The trace extension
theorem allows us to subtract a $1$-form $\beta\in H^s(M,T^*M)$ with
$\beta(\boldsymbol\nu)=\alpha(\boldsymbol\nu)$ and
$\|\beta\|_{H^s(M,T^*M)}\leq C\|\alpha(\boldsymbol\nu)\|_{H^{s-1/2}(\partial M)}$.
The difference, again denoted by $\alpha$, has zero normal component.
Integration by parts in the Weitzenböck identity controls its $H^1(M,T^*M)$ norm by
$d\alpha,\delta\alpha$ and its $L^2(M,T^*M)$ norm: the boundary term is a
quadratic form of order zero in the tangential trace. Its coefficients
are bounded, and the trace inequality
\[
 \|\operatorname{Tr}\alpha\|_{L^2(\partial M)}^2
 \leq\varepsilon\|\alpha\|_{H^1(M)}^2+
                     C_\varepsilon\|\alpha\|_{L^2(M)}^2
\]
allows it to be absorbed. Indeed, if the absolute value of the
boundary term is bounded by $C_\partial\|\operatorname{Tr}\alpha\|_{L^2(\partial M)}^2$,
it suffices to take $\varepsilon=[2(1+C_\partial)]^{-1}$.
In boundary normal coordinates, the first-order equations are
\[
 \partial_m\alpha_a=(d\alpha)_{ma}+\partial_a\alpha_m
                  \quad(1\leq a<m),
 \qquad
 \delta\alpha=-\partial_m\alpha_m
       -\sum_{a,b=1}^{m-1}g^{ab}\partial_a\alpha_b
       +\sum_{b=1}^m c_b\alpha_b,
\]
with smooth coefficients $c_b$. Once the tangential derivatives
have been estimated, these identities recover the normal derivatives.
Cutoffs, tangential estimates and a finite cover give
\eqref{eq:hodge-primer-orden-con-resto-L2}; its real-order Sobolev
versions are those of the cited regularity result.
When there is no boundary, one uses the interior elliptic estimate.

This regularity makes the vector fields in the kernel with zero normal
component smooth. The order-one estimate and Rellich's theorem show
that its $L^2(M,T^*M)$ unit ball is compact in $L^2(M,T^*M)$; thus $\mathcal H$ is
finite-dimensional. To replace the last term of
\eqref{eq:hodge-primer-orden-con-resto-L2}, suppose that
\eqref{eq:estimacion-primer-orden-hodge-fluidos} fails.
There would be vector fields $W_j$ with $\|W_j\|_{H^s(M,TM)}=1$ whose four
data on the right-hand side tend to zero. A subsequence converges in $L^2(M,TM)$.
Applying \eqref{eq:hodge-primer-orden-con-resto-L2} to its differences
makes it Cauchy in $H^s(M,TM)$. Its limit $W$ has $H^s(M,TM)$ norm one,
belongs to $\mathcal H$, and satisfies $\Pi_{\mathcal H}W=0$.
This contradiction proves the estimate.
\end{proof}

\begin{lemma}[A left inverse in fixed spaces]
\label{lem:inversa-izquierda-familia-operadores-fluidos}
Let $E,F$ be Hilbert spaces and let $D_0\in\mathcal L(E,F)$ satisfy
$\|D_0v\|_F\geq c\|v\|_E$, with $c>0$.
There is $L\in\mathcal L(F,E)$ such that $LD_0=I_E$.
If $z\mapsto D_z\in\mathcal L(E,F)$ is of class $C^q$, with
$q\in\mathbb N\cup\{\infty\}$ and $q\geq1$, then, in a neighborhood where
$\|L(D_z-D_0)\|_{\mathcal L(E)}<1/2$, every solution of
$D_zv=b_z$ satisfies
\begin{equation}
\label{eq:solucion-inversa-izquierda-fluidos}
 v=(LD_z)^{-1}Lb_z,
 \qquad \|(LD_z)^{-1}\|_{\mathcal L(E)}\leq2.
\end{equation}
In particular, it depends on the parameters with class $C^q$ if $b_z$ does.
The operator $D_z$ need not be surjective onto $F$.
\end{lemma}

\begin{proof}
The lower bound makes the image of $D_0$ closed. If $P_0$ is the
orthogonal projection of $F$ onto that image, take
$L=\bigl(D_0:E\to\operatorname{im}D_0\bigr)^{-1}P_0$; then
$\|L\|\leq c^{-1}$. Write $A_z=LD_z$.
The operator-norm convergent series
\[
 A_z^{-1}=\sum_{j=0}^{\infty}(I_E-A_z)^j
\]
gives the bound two. For $\|A_z^{-1}H\|<1/2$, the resolvent identity gives
\[
 \|(A_z+H)^{-1}-A_z^{-1}+A_z^{-1}HA_z^{-1}\|
 \leq2\|A_z^{-1}\|^3\|H\|^2.
\]
Thus inversion is differentiable, with derivative
$H\mapsto-A_z^{-1}HA_z^{-1}$. Iterating the product rule proves
its smoothness. Multiplying $D_zv=b_z$ by $L$ proves
\eqref{eq:solucion-inversa-izquierda-fluidos} and the asserted dependence.
\end{proof}

\begin{proposition}[Smoothness of the projection in Lagrangian coordinates]
\label{prop:proyeccion-lagrangiana-suave-compacta}
Under the preceding compact hypotheses, with $s>m/2+1$, the map
\[
 \mathbb P_\eta W
 :=\bigl(\mathbb P(W\circ\eta^{-1})\bigr)\circ\eta,
 \qquad
 W\in\mathcal E^s_\eta:=H^s(M,\eta^*TM),
 \quad \eta\in\operatorname{Diff}^s_\mu(M),
\]
is a smooth bundle morphism. It is the orthogonal projection, for the
$L^2(M,\eta^*TM)$ metric, of $\mathcal E^s_\eta$ onto
$T_\eta\operatorname{Diff}^s_\mu(M)$.
\end{proposition}

\begin{proof}
The projection is well defined and continuous in each fiber by
\eqref{eq:cotas-Q-P-hodge-leray} and the composition bounds.
The change of variables by the volume-preserving map $\eta$
proves orthogonality. It remains to establish dependence on $\eta$.

Choose an orthonormal basis $Z_1,\ldots,Z_b$ of $\mathcal H$ in
$L^2(M,TM)$. If $\mathcal H=\{0\}$, these data are omitted and
$\mathbb R^b=\{0\}$. For $W\in\mathcal E^s_\eta$, define
\[
 \begin{aligned}
 \mathcal C_\eta W&:=
       d\bigl((W\circ\eta^{-1})^\flat\bigr)\circ\eta,\\
 \mathcal B_\eta W&:=
       \bigl(\operatorname{div}_{\mathbf g}(W\circ\eta^{-1})\bigr)\circ\eta,\\
 \mathcal N_\eta W&:=
       \langle\operatorname{Tr}W,\boldsymbol\nu\circ\eta\rangle_{\mathbf g},\\
 \ell_{\eta,j}(W)&:=\int_M
       \langle W(x),Z_j(\eta(x))\rangle_{\mathbf g}\,d\lambda_{\mathbf g}(x),
                         \qquad 1\leq j\leq b.
 \end{aligned}
\]
The first composition means a section of the induced bundle
$\eta^*\Lambda^2T^*M$, without multiplication by $d\eta$.
The normal datum is defined only on $\partial M$ and is omitted
when there is no boundary. The combined operator is
\[
 \mathcal D_\eta W=(\mathcal C_\eta W,\mathcal B_\eta W,
                       \mathcal N_\eta W,\ell_\eta(W)),
\]
with values in
\[
 \mathcal F^{s-1}_\eta:=
 H^{s-1}(M,\eta^*\Lambda^2T^*M)\oplus H^{s-1}(M)
 \oplus H^{s-1/2}(\partial M)\oplus\mathbb R^b.
\]

Let us verify smoothness in operator norm. In coordinates
$y=\eta(x)$, write $\displaystyle W(x)=\sum_{a=1}^m W^a(x)\partial_{y^a}$
and $B(x)=(D\eta(x))^{-1}$. The components of the covariant derivative
of the Eulerian field, composed with $\eta$, are
\begin{equation}
\label{eq:derivada-conjugada-primer-orden-fluidos}
 A^a_b(\eta,W)
 =\sum_{j=1}^m B^j_b\partial_jW^a
      +\sum_{c=1}^m\Gamma^a_{bc}(\eta)W^c,
                    \qquad 1\leq a,b\leq m.
\end{equation}
Hence
\[
 \mathcal B_\eta W=\sum_{a=1}^m A^a_a(\eta,W),
 \qquad
 (\mathcal C_\eta W)_{ab}
   =\sum_{c=1}^m\bigl(g_{bc}(\eta)A^c_a(\eta,W)
                       -g_{ac}(\eta)A^c_b(\eta,W)\bigr).
\]
Only first derivatives of $\eta$ and $W$ occur.
Matrix inversion is smooth in a neighborhood where determinants
are bounded away from zero; $H^{s-1}(M)$ is an algebra. Superposition,
multiplication and the trace therefore give smoothness of the first
three components between the stated spaces. For $\ell_\eta$, one uses
smoothness of each $Z_j$, superposition in $H^s(M,\eta^*TM)$, and the continuous
$L^2(M,\eta^*TM)$ pairing. These assertions hold in operator norm because all
expressions are linear in $W$ and the product bounds control the
unit ball of $H^s(M,\eta^*TM)$ uniformly.

Work near the identity and trivialize the bundles by parallel transport
in a normal neighborhood of the diagonal. In the presence of boundary,
use a smooth extension of the metric to the double. These are smooth
trivializations of $\mathcal E^s$ and $\mathcal F^{s-1}$, since they
are constructed by superposition. Estimate
\eqref{eq:estimacion-primer-orden-hodge-fluidos} gives a left inverse
$L$ for $\mathcal D_{\operatorname{id}}$ in these fixed spaces.
Choose a neighborhood $\mathcal V$ where
\[
 \|L(\mathcal D_\eta-\mathcal D_{\operatorname{id}})\|<1/2.
\]
Operator-norm continuity guarantees its existence.

Let $V=\mathbb P_\eta W$. The Eulerian difference
$(W-V)\circ\eta^{-1}$ is a gradient. Its dual form is closed, and the field is
orthogonal to $\mathcal H$.
Moreover, $V$ has zero divergence and normal component.
Thus, exactly,
\[
 \mathcal D_\eta V=(\mathcal C_\eta W,0,0,\ell_\eta(W)).
\]
The preceding lemma, applied to this already existing solution, gives
\begin{equation}
\label{eq:proyeccion-lagrangiana-inversa-fija}
 \mathbb P_\eta W
   =(L\mathcal D_\eta)^{-1}L
                         (\mathcal C_\eta W,0,0,\ell_\eta(W)).
\end{equation}
Every factor is smooth in the chosen neighborhood. This proves smoothness
without differentiating inversion of diffeomorphisms at a fixed level.
Equivariance under right translation transports this construction from
the identity to every point of the group.
\end{proof}

\begin{theorem}[Kinetic connection and geodesic vector field]
\label{teo:spray-euler-compacto-demostrado}
For a smooth compact, connected and oriented manifold, with or without
boundary, and $s>m/2+1$, the kinetic metric on
$\operatorname{Diff}^s_\mu(M)$ has a smooth Levi--Civita connection
and a smooth spray. Every initial datum in its tangent bundle determines
a unique local geodesic on an interval $I=(-T,T)$, with $T>0$.
Its Eulerian velocity has the regularity stated
in Proposition~\ref{prop:euler-geodesica} and satisfies Euler, with
normalized pressure in $C(I,H^{s+1}(M))$.
\end{theorem}

\begin{proof}
In the space of maps into the double of $M$, or into $M$ when there
is no boundary, the group is a split submanifold. When a boundary is
present, the condition of mapping $\partial M$ to $\partial M$ is split
by the trace extension operator; within it, apply
Theorem~\ref{teo:ebin-marsden-diffmu-compacta}.
The pointwise connection $\overline\nabla$ is smooth. Its coordinate
terms are superpositions of the smooth Levi--Civita coefficients and
products of two $H^s(M,\eta^*TM)$ sections.
The connection on the subgroup is
\[
 \nabla^{\mu}_X Y:=\mathbb P_\eta(\overline\nabla_XY).
\]
The preceding proposition proves its smoothness. The pointwise connection
is compatible with the $L^2(M,\eta^*TM)$ integral. After projection, the orthogonal
components do not contribute to pairings with tangent vector fields,
so $\nabla^\mu$ is compatible with the kinetic metric.
Also,
\[
 \nabla^\mu_XY-\nabla^\mu_YX
       =\mathbb P_\eta[X,Y]=[X,Y].
\]
Its torsion vanishes. The Koszul identity determines a metric connection
with zero torsion uniquely, since the $L^2(M,\eta^*TM)$ pairing separates tangent
vectors.

In an ambient chart, let $x$ be the coordinate of $\eta$, let $v$
represent its velocity, and write $P_x$ for the projection in these
coordinates and $\Gamma_x$ for the ambient connection term.
The tangency condition is $P_xv=v$. Differentiating it and using
the equation $P_x(\dot v+\Gamma_x(v,v))=0$, we obtain
\begin{equation}
\label{eq:spray-euler-proyeccion-derivada}
 \mathbf S(x,v)=
       \bigl(v,DP_x[v]v-P_x\Gamma_x(v,v)\bigr).
\end{equation}
Indeed, differentiating $P_x^2=P_x$ gives
$P_xDP_x[v]v=0$ when $P_xv=v$; hence this formula satisfies
both identities and defines a vector tangent to the tangent bundle
of the subgroup. It is smooth and quadratic in $v$.
If $\pi:T\operatorname{Diff}^s_\mu(M)\to\operatorname{Diff}^s_\mu(M)$
is the projection, its first component expresses
$T\pi\circ\mathbf S=\operatorname{id}_{T\operatorname{Diff}^s_\mu(M)}$.

Let $E$ be the model space of a chart of this tangent bundle around
an initial datum $z_0$. Take $A=1+\|\mathbf S(z_0)\|$ and
$L=1+\|D\mathbf S(z_0)\|$. Continuity allows us to choose
a ball $B(z_0,\rho)$ contained in the domain such that
$\|\mathbf S(z)\|\leq A$ and $\|D\mathbf S(z)\|\leq L$ on that ball.
For $0<T<\min\{\rho/(2A),1/(2L)\}$, the operator
\[
 z\longmapsto\left(t\longmapsto
                    z_0+\int_0^t\mathbf S(z(\tau))\,d\tau\right)
\]
preserves the closed ball of radius $\rho/2$ around the constant
curve $z_0$ in $C([-T,T],E)$ and is a contraction.
Its fixed point is the unique local solution; this is the argument
of the ODE theorem in the chapter on Banach flows.
Identity \eqref{eq:aceleracion-euleriana-lagrangiana} transforms
$\nabla^\mu_t\dot\eta=0$ into
$\mathbb P(\partial_t\mathbf u+\nabla_{\mathbf u}\mathbf u)=0$.
Since $\partial_t\mathbf u$ is solenoidal, the Hodge decomposition
gives Euler with $p=-Q(\nabla_{\mathbf u}\mathbf u)$.
Proposition~\ref{prop:presion-ganancia-sobolev-compacta}
gives the higher regularity of the pressure.
\end{proof}

The connection and spray above are those of
\cite[Theorems~10.2 and~11.1, pp.~134--135; Appendix~A, pp.~153--157]{EbinMarsden1970}.
The first-order operator construction makes the dependence of the
projection on the Lagrangian chart explicit.
Propagation of additional regularity is studied in
Theorem~12.1 of \cite[p.~137]{EbinMarsden1970}.

The pressure is recovered by taking the divergence in
\eqref{eq:euler-incompresible}. The resulting equation is a Neumann
problem, whose boundary condition comes from the normal component of Euler's equation.
This elliptic problem is local in time, but local existence does not
imply completeness of the kinetic metric or global control of solutions.
In dimension three, global continuation of a smooth solution is a
question distinct from the geodesic interpretation.

On the complete manifold without boundary, let
$\Delta_1=d\delta+\delta d$ be the nonnegative Hodge Laplacian on
$L^2(M,T^*M)$, with its self-adjoint realization. Here $\delta$ is the formal
adjoint of $d$, and the complement of $\ker\Delta_1$ is taken in $L^2(M,T^*M)$.

\begin{theorem}[Eichhorn--Schmid: Euler in bounded geometry without boundary]
\label{teo:eichhorn-schmid-euler-abierto}
Let $(M,\mathbf g)$ (of dimension $m$) be smooth, oriented, connected, noncompact, and without
boundary, with bounded geometry of order $k$. Suppose that
\[
 \inf\sigma_{\mathrm e}\left(
 \Delta_1\restriction_{(\ker\Delta_1)^\perp}
 \right)>0,
 \qquad
 k-2\geq r>\frac{m}{2}+1,\qquad k,r\in\mathbb N.
\]
The identity component $\mathcal D_{\mu,0}^{2,r}(M)$ of the group of
volume-preserving Sobolev diffeomorphisms is a split
submanifold of class $C^{k+1-r}$ of $\mathcal D_0^{2,r}(M)$.
Each datum $\mathbf u_0\in H^r(M,TM)$ with
$\operatorname{div}_{\mathbf g}\mathbf u_0=0$ determines a unique local
Euler solution, with
\[
 \mathbf u\in C((-T,T),H^r(M,TM))
          \cap C^1((-T,T),H^{r-1}(M,TM)),
 \qquad \mathbf u(0)=\mathbf u_0,
\]
for some $T>0$. The pressure is chosen in $C((-T,T),H^{r+1}(M))$; its differential
therefore belongs to the closed range of $d:H^1(M)\to L^2(M,T^*M)$. This condition is part of the class of solutions in which
uniqueness is asserted. No mean normalization is imposed on the base of
infinite volume.
The kinetic metric has a connection and spray of class at least $C^1$
on this submanifold. If the geometry is bounded to all orders,
both are smooth. The Lagrangian flow is, respectively, the flow
of a $C^1$ or smooth ODE on the Hilbert tangent bundle.
\end{theorem}

\begin{proof}
The hypothesis $k-2\geq r$ allows us to use the diffeomorphism manifold from
Theorem~\ref{teo:eichhorn-difeomorfismos-abiertos}; in particular,
$k+1-r\geq3$. The submanifold result for the volume form is
\cite[Theorem~3.3]{EichhornSchmidForms1996}, applied to the identity
component. The spectral condition provides closed range for the
Hodge operators and a complement to the directions changing
volume. The norms and topology are those of the Sobolev charts at the
identity; there is no need to extend to all maps a uniformity with
radii independent of their bounds.

Proposition~\ref{prop:brecha-espectral-rango-inversa}, applied to
the restriction of $\Delta_1$ to the orthogonal complement of its kernel,
provides a positive gap and a continuous inverse there.
Let $\gamma>0$ be a gap constant for $\Delta_1$ on
$(\ker\Delta_1)^\perp$, and let
$\Delta_0=\delta d=-\Delta_{\mathbf g}$ in $L^2(M)$. If $v\in C_c^\infty(M)$,
the form $dv$ is orthogonal to the kernel of $\Delta_1$, and commutation of
$d$ with the Hodge Laplacians gives
\[
 \|\Delta_0v\|_{L^2(M)}^2
 =\langle\Delta_1dv,dv\rangle_{L^2(M,T^*M)}
 \geq\gamma\|dv\|_{L^2(M,T^*M)}^2
 =\gamma\langle\Delta_0v,v\rangle_{L^2(M)}.
\]
Theorem~\ref{teo:sobolev-fraccionario-y-nucleo-calor-variedad-riemanniana-completa-esencialmente-auto}
states that $C_c^\infty(M)$ is a core for $\Delta_0$.
Thus, for $v\in\mathcal D(\Delta_0)$, choose $v_j\in C_c^\infty(M)$
with $v_j\to v$ and $\Delta_0v_j\to\Delta_0v$ in $L^2(M)$.
The squared norms and inner product in the preceding inequality
converge, so it holds for $v\in\mathcal D(\Delta_0)$.
If $E_0$ is the spectral resolution of $\Delta_0$, this inequality is
\[
 \int_{[0,\infty)}\lambda(\lambda-\gamma)\,
                     d\langle E_0(\lambda)v,v\rangle_{L^2(M)}\geq0.
\]
A nonzero spectral projection on a compact interval
$[a,b]\subset(0,\gamma)$ would give a nonzero vector in
$\mathcal D(\Delta_0)$ for which the integral is strictly
negative. All such projections therefore vanish, and
$\sigma(\Delta_0)\subseteq\{0\}\cup[\gamma,\infty)$. Proposition~\ref{prop:laplaciano-abierto-volumen-nucleo-resolvente}
gives infinite volume and $\ker_{L^2}\Delta_0=\{0\}$ under the
bounded-geometry bounds in this statement. If zero belonged to the
spectrum, it would be isolated and hence an eigenvalue by
Theorem~\ref{teo:espectro-esencial-proyecciones-weyl}, a contradiction.
Thus $\sigma(\Delta_0)\subset[\gamma,\infty)$, and
Proposition~\ref{prop:brecha-espectral-rango-inversa} gives
\[
 \|\Delta_0^{-1}\|_{\mathcal L(L^2(M))}\leq\gamma^{-1}.
\]

The Hodge decomposition and uniform elliptic estimates give
\[
 \mathbb P\mathbf X
 =\mathbf X-\bigl(d\Delta_0^{-1}\delta(\mathbf X^\flat)\bigr)^\sharp.
\]
This projection is continuous on $H^{r-1}(M,TM)$ and $H^r(M,TM)$. For
$\mathbf X=\nabla_{\mathbf u}\mathbf u$, take
\[
 p=-\Delta_0^{-1}\delta\bigl((\nabla_{\mathbf u}\mathbf u)^\flat\bigr).
\]
The divergence cancellation computed above places the argument of
$\Delta_0^{-1}$ in $H^{r-1}(M)$; the elliptic estimate gives $p\in H^{r+1}(M)$.
The local bounds are summed over a uniform cover of finite
multiplicity, and the gap controls the norm of the inverse on $L^2(M)$. This choice
preserves orthogonality between the pressure and divergence-free vector fields.
Allowing the gradient of an arbitrary harmonic potential, without requiring it to
come from a function in $L^2(M)$, would change the class of solutions and is not the
pressure condition used in this theorem.

To justify Lagrangian regularity under these same geometric bounds,
consider the Hilbert space
\[
 \mathcal H=\{Z\in L^2(M,TM)\mid Z^\flat\in\ker\Delta_1\}
\]
with the $L^2(M,TM)$ norm, and its orthogonal projection $\Pi_{\mathcal H}$.
We do not assume that $\mathcal H$ is finite-dimensional.
The uniform estimates for the first-order Hodge system give
\[
 \|W\|_{H^r(M,TM)}\leq C_0\bigl(
   \|d(W^\flat)\|_{H^{r-1}(M,\Lambda^2T^*M)}+
   \|\operatorname{div}_{\mathbf g}W\|_{H^{r-1}(M)}+
   \|W\|_{L^2(M,TM)}\bigr).
\]
Indeed, on a normal ball the order-one estimate follows from the
Weitzenböck identity applied to a cutoff of $W$. On commuting
derivatives up to order $r-1$, the commutators contain curvature and
its derivatives up to order $r-2$, together with lower derivatives
of $W$. Induction controls these terms. Choose balls of a common
radius smaller than the injectivity radius, and cutoffs whose
derivatives up to order $r$ satisfy common bounds. The finite
multiplicity of the cover allows the squares of the local estimates
to be summed. The constant $C_0$ depends on the radius, multiplicity
and fixed geometric bounds, and not on the center of the ball.
The estimate extends by density to $H^r(M,TM)$ and has the corresponding
weak regularity version.

The spectral gap, applied to the quadratic form of $\Delta_1$, gives
\[
 \gamma\|(I-\Pi_{\mathcal H})W\|_{L^2(M,TM)}^2
 \leq\|d(W^\flat)\|_{L^2(M,\Lambda^2T^*M)}^2+
       \|\operatorname{div}_{\mathbf g}W\|_{L^2(M)}^2.
\]
Consequently,
\begin{equation}
\label{eq:control-primer-orden-uniforme-euler}
 \|W\|_{H^r(M,TM)}\leq C_1\bigl(
   \|d(W^\flat)\|_{H^{r-1}(M,\Lambda^2T^*M)}+
   \|\operatorname{div}_{\mathbf g}W\|_{H^{r-1}(M)}+
   \|\Pi_{\mathcal H}W\|_{L^2(M,TM)}\bigr).
\end{equation}
Weak regularity applied to $Z\in\mathcal H$ also gives
\[
 \|Z\|_{H^3(M,TM)}\leq C_{\mathcal H}\|Z\|_{L^2(M,TM)}.
\]
Here $d(Z^\flat)=0$ and $\operatorname{div}_{\mathbf g}Z=0$ weakly.
The required orders are available because $r\geq2$ and $k\geq r+2$.
With bounded geometry of all orders, the same estimate holds in
each $H^a(M,TM)$, with a constant depending on $a$.

For $\eta\in\mathcal D_{\mu,0}^{2,r}(M)$, define
\[
 J_\eta:\mathcal H\longrightarrow L^2(M,\eta^*TM),
 \qquad J_\eta Z=Z\circ\eta.
\]
The change of variables proves that it is an $L^2(M,TM)$ isometry.
Its Hilbert adjoint satisfies
\[
 J_\eta^*W=\Pi_{\mathcal H}(W\circ\eta^{-1}).
\]
Let us verify dependence of $J_\eta$ in operator norm without
differentiating this last composition. Trivialize the fibers by
parallel transport in a uniform neighborhood of the diagonal and
use a Sobolev chart of the identity, with parameter $z$ in the space
$E_\mu=\{h\in H^r(M,TM)\mid\operatorname{div}_{\mathbf g}h=0\}$.
On a
sufficiently small ball in the chart, the Jacobians are bounded
above and below, and the derivatives of its parametrizations up
to order three have local bounds. For $j=0,1,2,3$, the chain rule
and the parallel transport bounds give
\begin{equation}
\label{eq:derivadas-composicion-armonica-euler}
 \|D^jJ_z[h_1,\ldots,h_j]Z\|_{L^2(M,TM)}
 \leq C_j\|Z\|_{H^j(M,TM)}
                 \prod_{a=1}^{j}\|h_a\|_{H^r(M,TM)}.
\end{equation}
For $j=0$, the product is one. Each term contains a derivative of
$Z$ of order at most $j$, composed with $\eta$, multiplied by chart
directions or their derivatives with respect to the parameter.
The latter are bounded in $C^0$ by the Sobolev embedding;
the change of variables controls the factor involving $Z$ in $L^2(M,TM)$.
The formula is first obtained for smooth vector fields and extended
by density in $H^3(M,TM)$. Integrating the third derivative along
segments of the parameter ball gives, for some $C>0$,
\[
 \|D^2J_{z_1}-D^2J_{z_0}\|_{
   \mathcal L^2(E_\mu;\mathcal L(\mathcal H,L^2(M,TM)))}
 \leq C\|z_1-z_0\|_{H^r(M,TM)}.
\]
The directions $h_a$ belong to $E_\mu$, with the norm of $H^r(M,TM)$. The first- and second-order integral formulas prove that
$z\mapsto J_z$ is $C^2$ in $\mathcal L(\mathcal H,L^2(M,TM))$.
Taking adjoints is a continuous linear operation on these operator
spaces; it follows that $\eta\mapsto J_\eta^*$ is $C^2$ as an
operator from $H^r(M,\eta^*TM)$ to $\mathcal H$.

Now use $\mathcal C_\eta,\mathcal B_\eta$ defined in the compact
case, without the boundary datum, and set
\[
 \mathcal D_\eta W=(\mathcal C_\eta W,\mathcal B_\eta W,J_\eta^*W),
 \qquad
 \mathcal F^{r-1}_\eta=
 H^{r-1}(M,\eta^*\Lambda^2T^*M)\oplus H^{r-1}(M)\oplus\mathcal H.
\]
Formula \eqref{eq:derivada-conjugada-primer-orden-fluidos} shows that
the first two components are $C^2$ in operator norm. Let us specify
the orders involved: the expressions are estimated in $H^{r-1}(M,TM)$,
at most two parameter derivatives are taken, and the charts and
their trivializations require mixed jets of total order at most
$r+3$. Lemma~\ref{lem:jets-mixtos-exponencial-orden-finito} controls
them because $r+3\leq k+1$. The fixed coefficients, which need not
be integrable over $M$, act as uniformly bounded multipliers.
Their variations are estimated on the uniform cover using the
Sobolev product, since $r-1>m/2$. Summing the squares gives global
operator-norm bounds. This observation, together with the estimate
for $J_\eta^*$, gives $C^2$ dependence for $\mathcal D_\eta$.

Estimate \eqref{eq:control-primer-orden-uniforme-euler} and
Lemma~\ref{lem:inversa-izquierda-familia-operadores-fluidos} provide
a left inverse $L$ for $\mathcal D_{\operatorname{id}}$.
In a neighborhood where
$\|L(\mathcal D_\eta-\mathcal D_{\operatorname{id}})\|<1/2$, let
$V=\mathbb P_\eta W$, initially defined by conjugating the bounded
projection $\mathbb P$. The Eulerian difference is the gradient
of a function in $H^1(M)$; hence it is closed and orthogonal
to $\mathcal H$. Consequently,
\[
 \mathcal D_\eta V=(\mathcal C_\eta W,0,J_\eta^*W),
 \qquad
 \mathbb P_\eta W=(L\mathcal D_\eta)^{-1}
                         L(\mathcal C_\eta W,0,J_\eta^*W).
\]
The projection is $C^2$ near the identity, and right translations
extend this conclusion to the whole group. With geometric bounds
of all orders, the same argument applies to each derivative and
proves smoothness.

The pointwise connection followed by $\mathbb P_\eta$ is compatible
with the kinetic metric and has zero torsion, by the same calculation
as in Theorem~\ref{teo:spray-euler-compacto-demostrado}. The atlas is
at least $C^3$; the two derivatives of chart transitions involved
in the connection leave coefficients of class at least $C^1$.
Formula \eqref{eq:spray-euler-proyeccion-derivada} therefore defines
a spray of class at least $C^1$ under $k-2\geq r$. With bounded
geometry of all orders, the connection and spray are smooth.
The contraction argument of the same theorem gives a local
geodesic $\eta$ of class $C^2$ with
$\eta(0)=\operatorname{id}$ and $\dot\eta(0)=\mathbf u_0$.
Composition at lower levels gives
$\mathbf u=\dot\eta\circ\eta^{-1}\in C((-T,T);H^r(M,TM))\cap C^1((-T,T);H^{r-1}(M,TM))$,
on the interval of existence and over $M$.
The acceleration identity and the projection give Euler with the
pressure already constructed. The product bounds and elliptic
estimate applied to the difference of two solenoidal data give
\[
 \|p(\mathbf u)-p(\mathbf v)\|_{H^{r+1}(M)}
 \leq C\bigl(\|\mathbf u\|_{H^r(M,TM)}+
             \|\mathbf v\|_{H^r(M,TM)}\bigr)
             \|\mathbf u-\mathbf v\|_{H^r(M,TM)}.
\]
Thus $p\in C((-T,T),H^{r+1}(M))$.

Uniqueness also holds directly in the Eulerian class of the statement.
Let $\mathbf u,\mathbf v$ be two solutions and set
$\mathbf w=\mathbf u-\mathbf v$. Then
\[
 \partial_t\mathbf w+\nabla_{\mathbf u}\mathbf w+
       \nabla_{\mathbf w}\mathbf v=-\nabla(p(\mathbf u)-p(\mathbf v)).
\]
The pressure gradient is orthogonal to $\mathbf w$ in $L^2(M,TM)$.
Moreover, $|\mathbf w|_{\mathbf g}^2\in H^1(M)$ because
$\mathbf w\in H^r(M,TM)\subset C_b^1(M,TM)$. The vanishing weak
divergence of $\mathbf u$ and density of $C_c^\infty(M)$ in $H^1(M)$ give
\[
 \int_M\langle\nabla_{\mathbf u}\mathbf w,\mathbf w\rangle_{\mathbf g}
                  \,d\lambda_{\mathbf g}
 =\frac12\int_M\langle\mathbf u,
                  \nabla|\mathbf w|_{\mathbf g}^2\rangle_{\mathbf g}
                  \,d\lambda_{\mathbf g}=0.
\]
Since $\partial_t\mathbf w\in H^{r-1}(M,TM)\subset L^2(M,TM)$,
its squared norm can be differentiated, giving
\[
 \left|\frac{d}{dt}\|\mathbf w(t)\|_{L^2(M,TM)}^2\right|
 \leq2\|\nabla\mathbf v(t)\|_{L^\infty(M,T^*M\otimes TM)}
                          \|\mathbf w(t)\|_{L^2(M,TM)}^2.
\]
The coefficient is continuous and bounded on each compact time
interval. Gronwall, applied from the initial datum and also with
time reversed, implies $\mathbf w=0$ when $\mathbf w(0)=0$.
The pressures agree by uniqueness of $\Delta_0^{-1}$ in $L^2(M)$.
The resulting projection and flow are those of
\cite[Proposition~4.9, Theorem~4.3 and Corollary~4.6, pp.~172--175]{EichhornSchmidForms1996}.
\end{proof}

\section*{Loss of regularity and choice of variables}
\addcontentsline{toc}{section}{Loss of regularity and choice of variables}
\index{loss of derivatives}
\index{compactness!geometric obstructions}

The three examples show that a variational formulation does not by itself
determine the appropriate function space. In the harmonic map problem,
the tension field contains two derivatives of the map. In Yang--Mills theory, the curvature
loses one derivative and the gradient
\(d_{\mathbf{A}}^*\mathbf{F}_{\mathbf{A}}\) loses two. In Euler's equation, the
Eulerian expression loses one derivative, whereas Lagrangian
variables produce a smooth spray in the compact case. Under the
bounded geometry and spectral gap hypotheses of the preceding theorem,
the spray is at least $C^1$ with $k\geq r+2$, and smooth if the
geometric bounds hold at all orders.

Symmetries also determine the form of the analysis. Gauge invariance
requires separating the directions \(d_{\mathbf{A}}\boldsymbol{\xi}\); invariance under
reparametrizations leads to working on a diffeomorphism group; and
ellipticity controls local regularity, but does not by itself provide
global compactness. On a noncompact manifold, uniform conditions
and often spectral hypotheses enter the picture.

A strong metric may change the representative of a differential within
the chosen tangent space, but it does not eliminate the loss of regularity of the
original geometric operator. This distinction between the energy level, the
geometric level, and the level at which the equation lives will be essential in the
variational applications of the following chapters.

\part{Variational methods on Riemannian manifolds}
\chapter*{Introduction to Part VII}
\addcontentsline{toc}{chapter}{Introduction to Part VII}
\markboth{Part VII. Variational methods on Riemannian manifolds}{Introduction to Part VII}

Many nonlinear equations can be written as the condition that an energy remain unchanged to first order under perturbations of the unknown. In this case, solving the equation amounts to finding critical points of a functional defined on a space of functions or sections. Computing the first variation identifies the equation; the difficulty lies in proving that these critical points exist and in understanding the mechanisms that prevent the sequences used to find them from escaping.

The direct method studies minima. A minimizing sequence is controlled through coercivity, a subsequence is extracted by reflexivity, and the minimum is recovered through lower semicontinuity. Superlinear problems require a different geometry: the desired critical point may be a saddle point rather than a minimum. Palais--Smale conditions, pseudogradients, and deformation lemmas then allow minimax arguments to be constructed, including the mountain pass theorem.

The Nehari manifold introduces a third perspective. When the functional has a suitable radial geometry, each positive ray may intersect a natural constraint at a unique point. The problem is thus reduced to a more controlled geometry without artificially imposing an external condition. These chapters bring together the tools developed throughout the book: Sobolev spaces specify the domain of the functional, embeddings control the nonlinearities, integration by parts identifies the weak equation, and Banach geometry makes it possible to formulate constraints, tangent spaces, and deformation flows. Variational theory therefore emerges as a natural continuation of the analysis developed in the preceding parts.

\begin{semblanzaHistorica}{From Euler and Lagrange to Palais and Smale}
The calculus of variations began with problems about curves and surfaces in which a geometric or physical quantity had to be minimized. Euler and Lagrange developed the differential language that makes it possible to pass from an energy to an equation. With functional analysis, the unknown could be understood as a point of an infinite-dimensional space. The direct method organized existence arguments for minima through compactness and semicontinuity; Morse theory and minimax methods extended the search to critical points that are not minima; and Palais and Smale introduced a compactness condition adapted to these problems. The Nehari constraint, in turn, made a radial identity into a particularly effective tool for nonlinear equations.
\end{semblanzaHistorica}

\chapter{Variational formulation of elliptic problems}
\label{cap:formulacion-variacional-problemas-elipticos}
\glsadd{metodo-variacional}
\index{variational principle}

Variational calculus allows partial differential equations to be studied through ``energy'' functionals defined on spaces of functions or sections. This perspective applies both on Euclidean domains and on Riemannian manifolds and, more generally, on smooth vector bundles. The common idea is first to formulate a notion of weak solution, obtained by multiplying the equation by test functions or sections and integrating by parts, and then to construct an energy functional whose Fréchet derivative reproduces that weak identity exactly. In this way, weak solutions arise as critical points of the functional.

The Green and integration-by-parts formulas developed in the preceding chapters provide the link between a differential operator and its weak formulation. In this part we will use them to identify the test spaces and recognize the Fréchet derivative of the associated energy functional.

For Poisson's equation, the quadratic functional
\[
u\longmapsto
\frac{1}{2}\int_M|du|_{\mathbf{g}}^2\,d\lambda_{\mathbf{g}}
-\int_Mfu\,d\lambda_{\mathbf{g}}
\]
produces the weak equation $-\Delta_{\mathbf{g}}u=f$. Replacing the quadratic term
by $\displaystyle \frac{1}{p}\int_M|du|_{\mathbf{g}}^p\,d\lambda_{\mathbf{g}}$ gives the
$p$-Laplacian, and the natural space changes from $W^{1,2}(M)$ to $W^{1,p}(M)$.
Zeroth-order nonlinearities will be treated using Nemytski operators:
the Carathéodory conditions and growth bounds will determine
when the map $u\mapsto f(\,\cdot\,,u)$ acts continuously between the
required Lebesgue spaces.

\begin{semblanzaHistorica}{From an equation to an energy}
The calculus of variations arose from the study of curves that extremize lengths, times, or actions. In modern elliptic equations, the same principle is formulated on function spaces: the first variation of an energy produces a weak identity, and critical points represent solutions. This translation is powerful because it replaces a pointwise equation by global properties of the functional, but it works only when integrability, differentiability, and boundary terms are justified in the chosen space.
\end{semblanzaHistorica}

\section{Carathéodory functions and superposition operators}
\label{sec:caratheodory-variacional}

In a nonlinearity $f(x,\xi)$, the point $x\in M$ is the spatial variable and $\xi\in V$ represents the value of the unknown $u$ at that point. These two variables play different roles in variational problems. Often, the coefficients of the equation depend only measurably on the base point, whereas the dependence on the unknown must retain enough continuity for compositions to be measurable and integral functionals to be well defined. The Carathéodory conditions introduced below formalize precisely this situation.

\begin{definition}[Carathéodory function]
\label{def:funcion-caratheodory-variedad}\glsadd{funcion-caratheodory}
\index{Caratheodory function@Carathéodory function}
Let $(M,\mathbf{g})$ be a Riemannian manifold with or without boundary, and let $V$ and $W$ be finite-dimensional real vector spaces. A function $f\colon M\times V\longrightarrow W$ is a \textbf{Carathéodory function} if, for each $\xi\in V$, the function $x\longmapsto f(x,\xi)$ is $\lambda_{\mathbf{g}}$-measurable and, for almost every $x\in M$, the function $\xi\longmapsto f(x,\xi)$ is continuous.
\end{definition}

\begin{proposition}[Measurability of compositions]
\label{prop:medibilidad-composicion-caratheodory-variedad}
Let $(M,\mathbf{g})$ be a Riemannian manifold with or without boundary, let $V$ and $W$ be finite-dimensional real vector spaces, let
$f\colon M\times V\longrightarrow W$ be a Carathéodory function, and let
$u\colon M\longrightarrow V$ be a measurable function. Then the function
$x\longmapsto f(x,u(x))$ is $\lambda_{\mathbf{g}}$-measurable.
\end{proposition}

\begin{proof}
Since $V$ is finite-dimensional, there exists a sequence $(u_j)_{j\in\mathbb N}$ of measurable simple functions such that $u_j(x)\to u(x)$ for almost every $x\in M$. For each $j\in\mathbb N$, we can write
\[
u_j=\displaystyle\sum_{i=1}^{N_j}\xi_{i,j}\mathbf{1}_{E_{i,j}},
\]
where the sets $E_{1,j},\dots,E_{N_j,j}$ are measurable and pairwise disjoint, and $\xi_{i,j}\in V$ for $1\leq i\leq N_j$. Adding the set $M\setminus\displaystyle\bigcup_{i=1}^{N_j}E_{i,j}$ with value $0\in V$ if necessary, we may assume, without loss of generality, that $\displaystyle M=\bigcup_{i=1}^{N_j}E_{i,j}$. Therefore,
\[
f(x,u_j(x))
=
\displaystyle\sum_{i=1}^{N_j}
f(x,\xi_{i,j})\mathbf{1}_{E_{i,j}}(x).
\]
Since $x\longmapsto f(x,\xi_{i,j})$ is measurable for each $i$ and the sets $E_{1,j},\dots,E_{N_j,j}$ are measurable, the function $x\longmapsto f(x,u_j(x))$ is measurable.

We now know that there is a measurable set $N\subseteq M$, with $\lambda_{\mathbf{g}}(N)=0$, such that $f(x,\cdot)$ is continuous and $u_j(x)\to u(x)$ for every $x\in M\setminus N$. Consequently, $f(x,u_j(x))\to f(x,u(x))$ for every $x\in M\setminus N$. Since $\lambda_{\mathbf{g}}$ is complete, an almost everywhere pointwise limit of measurable functions with values in a finite-dimensional space is measurable. Therefore, $x\longmapsto f(x,u(x))$ is $\lambda_{\mathbf{g}}$-measurable.
\end{proof}

The preceding formulation extends naturally to nonlinearities acting on sections of vector bundles. In this case, both the unknown and its image belong to fibers that depend on the base point. The next definition expresses the Carathéodory condition without choosing trivializations.

\begin{definition}[Fiberwise Carathéodory map]
\label{def:mapeo-caratheodory-haces}
\index{Caratheodory map@Carathéodory map!fiberwise}
Let $(M,\mathbf{g})$ be a Riemannian manifold with or without boundary, and let
$\mathbf{F}\to M$ and $\mathbf{G}\to M$ be smooth real vector bundles of finite rank. A fiber-preserving map $\mathscr G\colon \mathbf{F}\longrightarrow \mathbf{G}$ is called a \textbf{fiberwise Carathéodory map} if the following conditions hold:
\begin{enumerate}
\item there exists a measurable set $M_0\subseteq M$, with $\lambda_{\mathbf{g}}(M\setminus M_0)=0$, such that $\mathscr G_x\colon \mathbf{F}_x\longrightarrow \mathbf{G}_x$ is continuous for every $x\in M_0$;
\item for every open set $U\subseteq M$ and every measurable section $\mathbf{u}\colon U\longrightarrow \mathbf{F}\restriction_U$, the section $x\longmapsto\mathscr G_x(\mathbf{u}(x))$ is measurable on $\mathbf{G}\restriction_U$.
\end{enumerate}
\end{definition}

The second condition in the preceding definition can be checked using trivializations. The next result shows that the intrinsic definition agrees with the usual Carathéodory condition for the local representations.

\begin{proposition}[Local characterization of Carathéodory maps]
\label{prop:caratheodory-haces-caracterizacion-local}
Let $(M,\mathbf{g})$ be a Riemannian manifold with or without boundary, let
$\mathbf{F}\to M$ and $\mathbf{G}\to M$ be smooth real vector bundles of finite
rank, and let $\mathscr G\colon \mathbf{F}\longrightarrow \mathbf{G}$ be a fiber-preserving map. Then $\mathscr G$ is a fiberwise Carathéodory map if and only if, for every open set $U\subseteq M$ over which $\mathbf{F}$ and $\mathbf{G}$ admit smooth trivializations
\[
\Phi_{\mathbf{F}}\colon \mathbf{F}\restriction_U\longrightarrow U\times\mathbb R^r,
\qquad
\Phi_{\mathbf{G}}\colon \mathbf{G}\restriction_U\longrightarrow U\times\mathbb R^s,
\]
the local representation $\mathscr G_U\colon U\times\mathbb R^r\longrightarrow\mathbb R^s$, determined by
\[
\Phi_{\mathbf{G}}\circ\mathscr G\circ\Phi_{\mathbf{F}}^{-1}(x,\xi)
=
\bigl(x,\mathscr G_U(x,\xi)\bigr),
\]
is a Carathéodory function.
\end{proposition}

\begin{proof}
Suppose first that $\mathscr G$ is a fiberwise Carathéodory map. Fix $\xi\in\mathbb R^r$ and consider the smooth section $\mathbf{u}_\xi\colon U\longrightarrow \mathbf{F}\restriction_U$ defined by $\mathbf{u}_\xi(x):=\Phi_{\mathbf{F}}^{-1}(x,\xi)$. By the second condition in Definition~\ref{def:mapeo-caratheodory-haces}, the section $x\longmapsto\mathscr G_x(\mathbf{u}_\xi(x))$ is measurable on $\mathbf{G}\restriction_U$. Its local representation with respect to $\Phi_{\mathbf{G}}$ is the function $x\longmapsto\mathscr G_U(x,\xi)$, which is therefore measurable.

On the other hand, there is a measurable set $M_0\subseteq M$, with $\lambda_{\mathbf{g}}(M\setminus M_0)=0$, such that $\mathscr G_x\colon \mathbf{F}_x\longrightarrow \mathbf{G}_x$ is continuous for every $x\in M_0$. Since the restrictions of $\Phi_{\mathbf{F}}$ and $\Phi_{\mathbf{G}}$ to the fibers are linear isomorphisms, the function $\xi\longmapsto\mathscr G_U(x,\xi)$ is continuous for every $x\in U\cap M_0$. Therefore, $\mathscr G_U$ is a Carathéodory function.

Conversely, suppose that the local representations of $\mathscr G$ are Carathéodory functions. Since $M$ is second countable, we can choose a countable cover $(U_k)_{k\in\mathbb N}$ by common trivialization domains for $\mathbf{F}$ and $\mathbf{G}$. For each $k\in\mathbb N$, there exists a measurable set $N_k\subseteq U_k$, with $\lambda_{\mathbf{g}}(N_k)=0$, such that the local representation $\mathscr G_{U_k}(x,\cdot)$ is continuous for every $x\in U_k\setminus N_k$. If $N:=\displaystyle\bigcup_{k\in\mathbb N}N_k$, then $\lambda_{\mathbf{g}}(N)=0$ and $\mathscr G_x\colon \mathbf{F}_x\longrightarrow \mathbf{G}_x$ is continuous for every $x\in M\setminus N$.

Now let $U\subseteq M$ be an open set, and let $\mathbf{u}\colon U\longrightarrow \mathbf{F}\restriction_U$ be a measurable section. For each $k\in\mathbb N$, let $u_k\colon U\cap U_k\longrightarrow\mathbb R^r$ be the local representation of $\mathbf{u}$. Since $\mathscr G_{U_k}$ is a Carathéodory function, Proposition~\ref{prop:medibilidad-composicion-caratheodory-variedad} shows that the function
\[
x\longmapsto
\mathscr G_{U_k}(x,u_k(x))
\]
is measurable on $U\cap U_k$. This function is the local representation of the section $x\longmapsto\mathscr G_x(\mathbf{u}(x))$. Since the sets $U\cap U_k$ form a countable cover of $U$, it follows that $x\longmapsto\mathscr G_x(\mathbf{u}(x))$ is a measurable section of $\mathbf{G}\restriction_U$. Therefore, $\mathscr G$ is a fiberwise Carathéodory map.
\end{proof}

\begin{remark}
To verify that a fiber-preserving map is fiberwise Carathéodory, it suffices to check that its local representations are Carathéodory functions with respect to a countable trivializing cover. Proposition~\ref{prop:caratheodory-haces-caracterizacion-local} also shows that this property does not depend on the chosen trivializations.

In particular, every continuous fiber-preserving map $\mathscr G\colon \mathbf{F}\longrightarrow \mathbf{G}$ is a fiberwise Carathéodory map. If $\mathbf{F}=M\times V$ and $\mathbf{G}=M\times W$ are trivial bundles, the definition reduces to the usual Carathéodory condition for a function $f\colon M\times V\longrightarrow W$.
\end{remark}

\begin{definition}[Nemytski operator between bundles]
\label{def:operador-nemytski-variacional}
\index{Nemytski operator}
\index{superposition operator@superposition operator!between bundles}
Let $(M,\mathbf{g})$ be a Riemannian manifold with or without boundary, let
$\mathbf{F}\to M$ and $\mathbf{G}\to M$ be smooth real vector bundles of finite
rank, and let $\mathscr G\colon \mathbf{F}\longrightarrow \mathbf{G}$ be a fiberwise Carathéodory map. The operator defined pointwise by
$\mathcal N_{\mathscr G}(\mathbf{u})(x):=\mathscr G_x(\mathbf{u}(x))$, whenever
the resulting section belongs to the target space under consideration, is
called the \textbf{Nemytski operator}, \textbf{superposition operator}, or
\textbf{substitution operator} associated with $\mathscr G$.
\end{definition}

The next result gives a growth condition ensuring that the superposition operator acts continuously between Lebesgue spaces of sections. The independent term belonging to $L^q(M)$ allows manifolds of infinite volume to be included.

\begin{theorem}[Continuity of superposition operators between bundles]
\label{teo:superposicion-fibra-a-fibra}
\index{superposition operator@superposition operator!continuity}
Let $(M,\mathbf{g})$ be a Riemannian manifold with or without boundary, let $\mathbf{F}\to M$ and $\mathbf{G}\to M$ be smooth real vector bundles of finite rank equipped with bundle metrics $\mathbf{h}_{\mathbf{F}}$ and $\mathbf{h}_{\mathbf{G}}$, and let $\mathscr G\colon \mathbf{F}\longrightarrow \mathbf{G}$ be a fiberwise Carathéodory map. Let $1\leq p,q<\infty$ and suppose that there exist a measurable function $a\colon M\longrightarrow[0,\infty)$, with $a\in L^q(M)$, and a constant $C\geq0$ such that
\[
|\mathscr G_x(\xi)|_{\mathbf{h}_{\mathbf{G}}}
\leq
a(x)+C|\xi|_{\mathbf{h}_{\mathbf{F}}}^{\frac{p}{q}}
\]
for almost every $x\in M$ and every $\xi\in \mathbf{F}_x$. Then the operator
\[
\mathcal N_{\mathscr G}\colon L^p(M,\mathbf{F})\longrightarrow L^q(M,\mathbf{G}),
\qquad
\mathcal N_{\mathscr G}(\mathbf{u})(x):=\mathscr G_x(\mathbf{u}(x)),
\]
is well defined, is continuous, and maps bounded sets to bounded sets.
\end{theorem}

\begin{proof}
Let $\mathbf{u}\in L^p(M,\mathbf{F})$. Definition~\ref{def:mapeo-caratheodory-haces} ensures that $\mathcal N_{\mathscr G}(\mathbf{u})$ is a measurable section of $\mathbf{G}$. Moreover, by the growth hypothesis and the inequality $(\alpha+\beta)^q\leq2^{q-1}(\alpha^q+\beta^q)$, we have
\[
|\mathscr G_x(\mathbf{u}(x))|_{\mathbf{h}_{\mathbf{G}}}^q
\leq
2^{q-1}a(x)^q
+
2^{q-1}C^q|\mathbf{u}(x)|_{\mathbf{h}_{\mathbf{F}}}^p
\]
for almost every $x\in M$. Therefore,
\[
\|\mathcal N_{\mathscr G}(\mathbf{u})\|_{L^q(M,\mathbf{G})}^q
\leq
2^{q-1}\|a\|_{L^q(M)}^q
+
2^{q-1}C^q\|\mathbf{u}\|_{L^p(M,\mathbf{F})}^p.
\]
In particular, $\mathcal N_{\mathscr G}(\mathbf{u})\in L^q(M,\mathbf{G})$. Furthermore, if two sections agree almost everywhere, their images under $\mathcal N_{\mathscr G}$ also agree almost everywhere, so the operator is well defined on the equivalence classes of $L^p(M,\mathbf{F})$.

The preceding estimate also shows that $\mathcal N_{\mathscr G}$ maps bounded sets to bounded sets. If $B\subseteq L^p(M,\mathbf{F})$ is bounded, there exists $R>0$ such that $\|\mathbf{u}\|_{L^p(M,\mathbf{F})}\leq R$ for every $\mathbf{u}\in B$, and then
\[
\|\mathcal N_{\mathscr G}(\mathbf{u})\|_{L^q(M,\mathbf{G})}^q
\leq
2^{q-1}\|a\|_{L^q(M)}^q
+
2^{q-1}C^qR^p
\]
for every $\mathbf{u}\in B$.

We prove continuity. Let $(\mathbf{u}_j)_{j\in\mathbb N}\subseteq L^p(M,\mathbf{F})$ be such that $\mathbf{u}_j\to \mathbf{u}$ in $L^p(M,\mathbf{F})$. Consider an arbitrary subsequence $(\mathbf{u}_{j_\ell})_{\ell\in\mathbb N}$. Passing to a further subsequence, denoted by $(\mathbf{v}_\ell)_{\ell\in\mathbb N}$, we may assume that $\|\mathbf{v}_\ell-\mathbf{u}\|_{L^p(M,\mathbf{F})}\leq2^{-\ell}$ for every $\ell\in\mathbb N$.

Define
\[
H(x):=
\displaystyle\sum_{\ell=1}^{\infty}
|\mathbf{v}_\ell(x)-\mathbf{u}(x)|_{\mathbf{h}_{\mathbf{F}}}.
\]
Minkowski's inequality (Theorem~\ref{teo:b6-espacios-lp-espacios-de-lebesgue}) implies
\[
\|H\|_{L^p(M)}
\leq
\displaystyle\sum_{\ell=1}^{\infty}
\|\mathbf{v}_\ell-\mathbf{u}\|_{L^p(M,\mathbf{F})}
\leq
\displaystyle\sum_{\ell=1}^{\infty}2^{-\ell}
<
\infty.
\]
Consequently, there exists a measurable set $N\subseteq M$, with
$\lambda_{\mathbf{g}}(N)=0$, such that $H(x)<\infty$ for every
$x\in M\setminus N$. At each such point, convergence of the series
defining $H(x)$ implies that
$|\mathbf{v}_\ell(x)-\mathbf{u}(x)|_{\mathbf{h}_{\mathbf{F}}}\to0$.
Therefore, $\mathbf{v}_\ell(x)\to\mathbf{u}(x)$ in $\mathbf{F}_x$ and
\[
|\mathbf{v}_\ell(x)|_{\mathbf{h}_{\mathbf{F}}}
\leq
|\mathbf{u}(x)|_{\mathbf{h}_{\mathbf{F}}}+H(x).
\]
Continuity of $\mathscr G_x$ for almost every $x\in M$ gives
\[
\mathscr G_x(\mathbf{v}_\ell(x))
\longrightarrow
\mathscr G_x(\mathbf{u}(x))
\]
for almost every $x\in M$.

On the other hand, the triangle inequality and the growth hypothesis imply the existence of a constant $C_0>0$, independent of $\ell$, such that
\[
|\mathscr G_x(\mathbf{v}_\ell(x))-\mathscr G_x(\mathbf{u}(x))|_{\mathbf{h}_{\mathbf{G}}}^q
\leq
C_0\bigl(
a(x)^q
+
|\mathbf{v}_\ell(x)|_{\mathbf{h}_{\mathbf{F}}}^p
+
|\mathbf{u}(x)|_{\mathbf{h}_{\mathbf{F}}}^p
\bigr)
\]
for almost every $x\in M$. Since
\[
|\mathbf{v}_\ell(x)|_{\mathbf{h}_{\mathbf{F}}}^p
\leq
2^{p-1}
\bigl(
|\mathbf{u}(x)|_{\mathbf{h}_{\mathbf{F}}}^p+H(x)^p
\bigr),
\]
there is a constant $C_1>0$, independent of $\ell$, such that
\[
|\mathscr G_x(\mathbf{v}_\ell(x))-\mathscr G_x(\mathbf{u}(x))|_{\mathbf{h}_{\mathbf{G}}}^q
\leq
C_1\bigl(
a(x)^q
+
|\mathbf{u}(x)|_{\mathbf{h}_{\mathbf{F}}}^p
+
H(x)^p
\bigr)
\]
for almost every $x\in M$. The right-hand side belongs to $L^1(M)$. Theorem~\ref{convergencia dominada} gives
\[
\|\mathcal N_{\mathscr G}(\mathbf{v}_\ell)
-
\mathcal N_{\mathscr G}(\mathbf{u})\|_{L^q(M,\mathbf{G})}^q
\longrightarrow0.
\]

We have shown that every subsequence of $(\mathbf{u}_j)_{j\in\mathbb N}$ has a subsequence whose image converges to $\mathcal N_{\mathscr G}(\mathbf{u})$ in $L^q(M,\mathbf{G})$. By the subsequence criterion, $\mathcal N_{\mathscr G}(\mathbf{u}_j)\to\mathcal N_{\mathscr G}(\mathbf{u})$ in $L^q(M,\mathbf{G})$. Therefore, $\mathcal N_{\mathscr G}$ is continuous.
\end{proof}

The case of trivial bundles yields the usual result for vector-valued Carathéodory functions.

\begin{corollary}[Nemytski operator between Lebesgue spaces]
\label{cor:nemytski-crecimiento-variacional}
Let $(M,\mathbf{g})$ be a Riemannian manifold with or without boundary, let $1\leq p,q<\infty$, and let $f\colon M\times\mathbb R^N\longrightarrow\mathbb R^L$ be a Carathéodory function. Suppose that there exist a measurable function $a\colon M\longrightarrow[0,\infty)$, with $a\in L^q(M)$, and a constant $C\geq0$ such that
\[
\|f(x,\xi)\|
\leq
a(x)+C\|\xi\|^{\frac{p}{q}}
\]
for almost every $x\in M$ and every $\xi\in\mathbb R^N$. Then the operator
\[
\mathcal N_f\colon L^p(M,\mathbb R^N)
\longrightarrow
L^q(M,\mathbb R^L),
\qquad
\mathcal N_f(u)(x):=f(x,u(x)),
\]
is well defined, is continuous, and maps bounded sets to bounded sets.
\end{corollary}

\begin{proof}
Simply apply Theorem~\ref{teo:superposicion-fibra-a-fibra} to the trivial
bundles $\mathbf{F}=M\times\mathbb R^N$ and
$\mathbf{G}=M\times\mathbb R^L$, equipped with their Euclidean bundle
metrics.
\end{proof}

\begin{corollary}[The case of conjugate exponents]
\label{cor:superposicion-exponentes-conjugados-haces}
Let $(M,\mathbf{g})$ be a compact Riemannian manifold with or without boundary, let $\mathbf{F}\to M$ and $\mathbf{G}\to M$ be smooth real vector bundles of finite rank equipped with bundle metrics $\mathbf{h}_{\mathbf{F}}$ and $\mathbf{h}_{\mathbf{G}}$, and let $\mathscr G\colon \mathbf{F}\longrightarrow \mathbf{G}$ be a continuous fiber-preserving map. Suppose that $m>1$ and that there exists a constant $C>0$ such that
\[
|\mathscr G_x(\xi)|_{\mathbf{h}_{\mathbf{G}}}
\leq
C\bigl(1+|\xi|_{\mathbf{h}_{\mathbf{F}}}^{m-1}\bigr)
\]
for every $x\in M$ and every $\xi\in \mathbf{F}_x$. If $m':=\displaystyle\frac{m}{m-1}$, then the operator $\mathcal N_{\mathscr G}\colon L^m(M,\mathbf{F})\longrightarrow L^{m'}(M,\mathbf{G})$ is well defined, is continuous, and maps bounded sets to bounded sets.
\end{corollary}

\begin{proof}
Since $\mathscr G$ is continuous and preserves fibers, it is a fiberwise Carathéodory map. Moreover, $\displaystyle\frac{m}{m'}=m-1$. Since $M$ is compact, the constant function $a(x):=C$ belongs to $L^{m'}(M)$. The result follows from Theorem~\ref{teo:superposicion-fibra-a-fibra} by taking $p=m$ and $q=m'$.
\end{proof}

\begin{remark}[The compact case]\label{rem: nemytski-variedad-compacta-haz}
If $M$ is compact, then $\lambda_{\mathbf{g}}(M)<\infty$ and every constant function belongs to $L^q(M)$. Consequently, in Theorem~\ref{teo:superposicion-fibra-a-fibra} the growth hypothesis can be replaced by
\[
|\mathscr G_x(\xi)|_{\mathbf{h}_{\mathbf{G}}}
\leq
C\bigl(1+|\xi|_{\mathbf{h}_{\mathbf{F}}}^{\frac{p}{q}}\bigr).
\]
Similarly, in Corollary~\ref{cor:nemytski-crecimiento-variacional} it suffices to assume that $\|f(x,\xi)\|\leq C(1+\|\xi\|^{\frac{p}{q}})$ for almost every $x\in M$ and every $\xi\in\mathbb R^N$.
\end{remark}
\begin{corollary}[Differentiability of powers of the norm]
\label{cor:potencia-norma-Lp-C1-haces}
Let $(M,\mathbf{g})$ be a Riemannian manifold with or without boundary, and let $\mathbf{E}\to M$ be a smooth real vector bundle of finite rank, equipped with a bundle metric $\mathbf{h}_{\mathbf{E}}$. For $1<\ell<\infty$, the functional
\[
 Q_\ell(\mathbf{u}):=\frac1\ell\|\mathbf{u}\|_{L^\ell(M,\mathbf{E})}^\ell
\]
is of class $C^1$ on $L^\ell(M,\mathbf{E})$ and
\[
 DQ_\ell(\mathbf{u})(\mathbf{v})
 =\int_M |\mathbf{u}|_{\mathbf{h}_{\mathbf{E}}}^{\ell-2}
 \langle\mathbf{u},\mathbf{v}\rangle_{\mathbf{h}_{\mathbf{E}}}\,d\lambda_{\mathbf{g}}.
\]
The integrand is interpreted as zero where $\mathbf{u}=0$.
\end{corollary}

\begin{proof}
In each fiber, the map $\xi\mapsto |\xi|^{\ell}/\ell$ is of class $C^1$ and its differential is represented by $b_\ell(\xi):=|\xi|^{\ell-2}\xi$, with $b_\ell(0)=0$. At zero, differentiability follows from $|\xi|^\ell/|\xi|\to0$, since $\ell>1$; away from zero, apply the chain rule. Theorem~\ref{teo:superposicion-fibra-a-fibra}, with exponents $\ell$ and $\ell':=\ell/(\ell-1)$ and zero independent term, shows that $b_\ell$ is continuous from $L^\ell(M,\mathbf{E})$ to $L^{\ell'}(M,\mathbf{E})$.

By Hölder, the proposed formula defines a continuous linear functional of norm at most $\|\mathbf{u}\|_{L^\ell(M,\mathbf E)}^{\ell-1}$. The fundamental theorem of calculus in each fiber and Fubini give
\[
 Q_\ell(\mathbf{u}+\mathbf{h})-Q_\ell(\mathbf{u})
 =\int_0^1\int_M\langle b_\ell(\mathbf{u}+t\mathbf{h}),\mathbf{h}\rangle_{\mathbf{h}_{\mathbf{E}}}\,d\lambda_{\mathbf{g}}\,dt.
\]
Indeed, the integral of the absolute value is at most
$(\|\mathbf{u}\|_{L^\ell(M,\mathbf E)}+\|\mathbf{h}\|_{L^\ell(M,\mathbf E)})^{\ell-1}\|\mathbf{h}\|_{L^\ell(M,\mathbf E)}$, by Hölder. Subtracting the proposed functional, we obtain
\[
 \begin{aligned}
 &\left|Q_\ell(\mathbf{u}+\mathbf{h})-Q_\ell(\mathbf{u})
 -\int_M\langle b_\ell(\mathbf{u}),\mathbf{h}\rangle_{\mathbf{h}_{\mathbf{E}}}\,d\lambda_{\mathbf{g}}\right|\\
 &\qquad\leq\|\mathbf{h}\|_{L^\ell(M,\mathbf E)}
 \sup_{t\in[0,1]}\|b_\ell(\mathbf{u}+t\mathbf{h})-b_\ell(\mathbf{u})\|_{L^{\ell'}(M,\mathbf E)}.
 \end{aligned}
\]
Continuity of $b_\ell$ at $\mathbf{u}$ makes the supremum tend to zero as $\mathbf{h}\to0$: all the points $\mathbf{u}+t\mathbf{h}$ lie at distance at most $\|\mathbf{h}\|_{L^\ell(M,\mathbf E)}$ from $\mathbf{u}$. This proves Fréchet differentiability. Finally, Hölder gives
\[
 \|DQ_\ell(\mathbf{u})-DQ_\ell(\mathbf{w})\|_{(L^\ell)'}
 \leq\|b_\ell(\mathbf{u})-b_\ell(\mathbf{w})\|_{L^{\ell'}(M,\mathbf E)},
\]
so the differential is continuous.
\end{proof}

\section{Classical, strong, and weak solutions}
\label{sec:nociones-solucion-problema-modelo}

First consider the scalar problem
\begin{equation}\label{eq:problema-modelo-dirichlet}
\begin{cases}
-\Delta_{\mathbf{g}} u=f(x,u) & \text{in }M,\\
u=0 & \text{on }\partial M,
\end{cases}
\end{equation}
where $(M,\mathbf{g})$ is a compact Riemannian manifold of dimension $n$ with nonempty smooth boundary and $f\colon M\times\mathbb R\longrightarrow\mathbb R$ is a Carathéodory function. This equation will serve as a model for introducing the basic notions. We adopt the convention
\[
2^*:=
\begin{cases}
\dfrac{2n}{n-2},&n\geq3,\\[2mm]
+\infty,&n\in\{1,2\}.
\end{cases}
\]
When $|f(x,t)|$ grows like $|t|^{r-1}$, the growth is said to be sublinear if $1<r<2$ and superlinear subcritical if $2<r<2^*$. If $n\geq3$, the case $r=2^*$ is called critical and the case $r>2^*$ supercritical. For $n\in\{1,2\}$, no finite critical exponent occurs in the first-order embeddings. These terms depend on the function space under consideration: $2^*$ is the limiting exponent for the embedding of $W_0^{1,2}(M)$.

\begin{definition}[Classical solution]
Let $(M,\mathbf{g})$ be a compact Riemannian manifold with nonempty boundary,
and let $f\colon M\times\mathbb R\to\mathbb R$ be a Carathéodory function. A
classical solution of \eqref{eq:problema-modelo-dirichlet} is a function
$u\in C^2(M)$ satisfying the equation at every point of $M$ and whose
restriction to $\partial M$ is zero.
\end{definition}

\begin{definition}[Strong solution]
Let $(M,\mathbf{g})$ be a compact Riemannian manifold with nonempty boundary,
and let $f\colon M\times\mathbb R\to\mathbb R$ be a Carathéodory function. A
strong solution of \eqref{eq:problema-modelo-dirichlet} is a function
$u\in W^{2,2}(M)\cap W_0^{1,2}(M)$ such that
$-\Delta_{\mathbf{g}}u=f(x,u)$ almost everywhere on $M$.
\end{definition}

Green's formula, obtained earlier from Theorem~\ref{teo:divergencia}, shows that every sufficiently regular solution satisfies
\[
\int_M\langle\nabla u,\nabla v\rangle_{\mathbf{g}}\,d\lambda_{\mathbf{g}}
=
\int_M f(x,u)v\,d\lambda_{\mathbf{g}}.
\]
This leads to the notion natural in Sobolev spaces.

\begin{definition}[Weak solution]
\label{def:solucion-debil-problema-modelo}\glsadd{solucion-debil}
Let $(M,\mathbf{g})$ be a compact Riemannian manifold with nonempty boundary,
and let $f\colon M\times\mathbb R\to\mathbb R$ be a Carathéodory function. A
weak solution of \eqref{eq:problema-modelo-dirichlet} is a function
$u\in W_0^{1,2}(M)$ such that, for every $v\in W_0^{1,2}(M)$,
the function $x\mapsto f(x,u(x))v(x)$ belongs to $L^1(M)$ and
\[
\int_M\langle\nabla u,\nabla v\rangle_{\mathbf{g}}\,d\lambda_{\mathbf{g}}
=
\int_M f(x,u)v\,d\lambda_{\mathbf{g}}.
\]
\end{definition}

Every classical solution is both strong and weak. Passing from a weak solution to a strong or classical solution requires elliptic regularity results and depends on the geometry of $M$, the regularity of the boundary, and the hypotheses on the nonlinearity.

The same procedure applies to sections. Let $\mathbf{E}\to M$ be a smooth real vector bundle of finite rank, equipped with a bundle metric $\mathbf{h}_{\mathbf{E}}$ and a connection $\nabla^{\mathbf{E}}$ compatible with $\mathbf{h}_{\mathbf{E}}$. For a fixed $1\leq q<\infty$, one considers test sections
$\boldsymbol{\phi}\in W_0^{1,q}(M,\mathbf{E})$, and integration by parts
is performed using the induced connection and its formal adjoint.
The boundary condition is incorporated through the characterization of
$W_0^{1,q}(M,\mathbf{E})$ as the kernel of the trace operator,
given by Theorem~\ref{teo:caracterizacion-W0-trazas} with $m=1$.

\section{Energy functionals and critical points}
\label{sec:funcional-energia-puntos-criticos}
\glsadd{funcional-energia}\glsadd{punto-critico}

Let $X$ be a real Banach space, and let $I\colon X\longrightarrow\mathbb R$. We use the notation $DI(u)\in X'$ for the Fréchet derivative. A point $u\in X$ is critical if $DI(u)=0$.

\begin{proposition}[Scalar energy functional]
\label{prop:funcional-energia-escalar-C1}
Let $(M,\mathbf{g})$ be a compact Riemannian manifold of dimension $n$ with
nonempty boundary. Suppose that
$f\colon M\times\mathbb R\longrightarrow\mathbb R$ is a Carathéodory
function and that there exist $C>0$ and $r\in(1,2^*)$ such that
$|f(x,t)|\leq C(1+|t|^{r-1})$. Let $M_0\subseteq M$ be a measurable
set such that $\lambda_{\mathbf{g}}(M\setminus M_0)=0$ and
$f(x,\cdot)$ is continuous for every $x\in M_0$. Define $F(x,t):=\displaystyle\int_0^t f(x,s)\,ds$
for $x\in M_0$ and $F(x,t):=0$ for $x\in M\setminus M_0$, and set
\[
I(u):=\frac{1}{2}\int_M|\nabla u|_{\mathbf{g}}^2\,d\lambda_{\mathbf{g}}-\int_MF(x,u)\,d\lambda_{\mathbf{g}}
\]
on $X=W_0^{1,2}(M)$. Then $I\in C^1(X,\mathbb R)$ and
\[
DI(u)(v)=\int_M\langle\nabla u,\nabla v\rangle_{\mathbf{g}}\,d\lambda_{\mathbf{g}}-\int_Mf(x,u)v\,d\lambda_{\mathbf{g}}
\]
for all $u,v\in X$. In particular, the critical points of $I$ are exactly the weak solutions of \eqref{eq:problema-modelo-dirichlet}.
\end{proposition}

\begin{proof}
Write $I=Q-P$, where
\[
Q(u):=\frac{1}{2}\int_M|\nabla u|_{\mathbf{g}}^2\,d\lambda_{\mathbf{g}},
\qquad
P(u):=\int_MF(x,u)\,d\lambda_{\mathbf{g}}.
\]
The functional $Q$ is of class $C^1$ and
\[
DQ(u)(v)=\int_M\langle\nabla u,\nabla v\rangle_{\mathbf{g}}\,d\lambda_{\mathbf{g}}.
\]

The growth hypothesis implies
\[
|F(x,t)|
\leq
C|t|+\frac{C}{r}|t|^r.
\]
Since $X\hookrightarrow L^1(M)\cap L^r(M)$ continuously, $P$ is well defined. Likewise, Corollary~\ref{cor:nemytski-crecimiento-variacional}, applied with $p=r$ and $q=r':=\displaystyle\frac{r}{r-1}$, shows that
\[
u\longmapsto f(\cdot,u)
\]
is continuous from $L^r(M)$ to $L^{r'}(M)$.

For $u\in X$, define $K_u\in X'$ by
\[
K_u(v):=\int_Mf(x,u)v\,d\lambda_{\mathbf{g}}.
\]
If $C_S>0$ is a constant for the embedding $X\hookrightarrow L^r(M)$, then
\[
|K_u(v)|
\leq
\|f(\cdot,u)\|_{L^{r'}(M)}\|v\|_{L^r(M)}
\leq
C_S\|f(\cdot,u)\|_{L^{r'}(M)}\|v\|_X,
\]
so $K_u$ is linear and continuous.

Let $h\in X$. For almost every $x\in M$, Theorem~\ref{teo:b5-fundamental-calculo-riemann-banach}, applied to $s\longmapsto F(x,u(x)+sh(x))$, gives
\[
F(x,u+h)-F(x,u)
=
\int_0^1f(x,u+sh)h\,ds.
\]
Hölder's inequality (Proposition~\ref{desigualdad de holder}) and the growth bound justify applying Fubini's theorem for Bochner integrals~\ref{teo:b5-fubini-bochner}, in its scalar case, to $(M,\lambda_{\mathbf{g}})\times([0,1],ds)$. Consequently,
\[
P(u+h)-P(u)-K_u(h)
=
\int_0^1\int_M\bigl(f(x,u+sh)-f(x,u)\bigr)h\,d\lambda_{\mathbf{g}}\,ds
\]
and
\[
|P(u+h)-P(u)-K_u(h)|
\leq
C_S\|h\|_X
\int_0^1
\|f(\cdot,u+sh)-f(\cdot,u)\|_{L^{r'}(M)}\,ds.
\]
Let $(h_j)_{j\in\mathbb N}\subseteq X$ be such that $h_j\to0$ in $X$. For each $s\in[0,1]$, we have $u+sh_j\to u$ in $L^r(M)$, and continuity of the Nemytski operator implies
\[
\|f(\cdot,u+sh_j)-f(\cdot,u)\|_{L^{r'}(M)}\longrightarrow0.
\]
Moreover, there exists $R>0$ such that $\|u+sh_j\|_{L^r(M)}\leq R$ for every $j$ and every $s\in[0,1]$. The estimate obtained in the first part of the proof of Corollary~\ref{cor:nemytski-crecimiento-variacional} provides a constant $C_R>0$ such that
\[
\|f(\cdot,u+sh_j)-f(\cdot,u)\|_{L^{r'}(M)}\leq C_R
\]
for every $j$ and every $s\in[0,1]$. Theorem~\ref{convergencia dominada}, applied in the variable $s$, gives
\[
\int_0^1
\|f(\cdot,u+sh_j)-f(\cdot,u)\|_{L^{r'}(M)}\,ds
\longrightarrow0.
\]
Therefore,
\[
\lim_{h\to0}
\frac{|P(u+h)-P(u)-K_u(h)|}{\|h\|_X}=0,
\]
so $P$ is Fréchet differentiable and $DP(u)=K_u$.

Finally, if $u_j\to u$ in $X$, then $u_j\to u$ in $L^r(M)$ and
\[
\|DP(u_j)-DP(u)\|_{X'}
\leq
C_S\|f(\cdot,u_j)-f(\cdot,u)\|_{L^{r'}(M)}
\longrightarrow0.
\]
Thus, $P\in C^1(X,\mathbb R)$, $I=Q-P\in C^1(X,\mathbb R)$, and the asserted formula for $DI$ holds. Comparing it with Definition~\ref{def:solucion-debil-problema-modelo}, we conclude that $DI(u)=0$ if and only if $u$ is a weak solution of \eqref{eq:problema-modelo-dirichlet}.
\end{proof}

\begin{proposition}[Second derivative of the potential functional]
\label{prop:funcional-energia-escalar-C2}
Let $(M,\mathbf{g})$ be a compact Riemannian manifold of dimension $n$ with
nonempty boundary, and let $f$ be the nonlinearity in
Proposition~\ref{prop:funcional-energia-escalar-C1}. Suppose in addition that
$f(x,\cdot)$ is of class $C^1$ for almost every $x$, that $f_t$ is a Carathéodory
function, and that there exist $C>0$ and $r\in(2,2^*)$ such that
\[
|f_t(x,t)|
\leq
C\bigl(1+|t|^{r-2}\bigr)
\]
for almost every $x\in M$ and every $t\in\mathbb R$. Then $I\in C^2(X,\mathbb R)$ and
\[
D^2I(u)(v,w)
=
\int_M\langle\nabla v,\nabla w\rangle_{\mathbf{g}}\,d\lambda_{\mathbf{g}}
-
\int_M f_t(x,u)vw\,d\lambda_{\mathbf{g}}.
\]
\end{proposition}

\begin{proof}
Let $P(u):=\displaystyle\int_MF(x,u)\,d\lambda_{\mathbf{g}}$. By the preceding proposition,
\[
DP(u)(v)=\int_Mf(x,u)v\,d\lambda_{\mathbf{g}}.
\]
Define $s:=\displaystyle\frac{r}{r-2}$. Then $\frac{1}{s}+\frac{2}{r}=1$. Corollary~\ref{cor:nemytski-crecimiento-variacional}, applied with source exponent $r$ and target exponent $s$, shows that
\[
u\longmapsto f_t(\cdot,u)
\]
is continuous from $L^r(M)$ to $L^s(M)$, since $\frac{r}{s}=r-2$.

For $u\in X$, define the bilinear form
\[
B_u(h,v):=\int_Mf_t(x,u)hv\,d\lambda_{\mathbf{g}}.
\]
By Hölder's inequality (Proposition~\ref{desigualdad de holder}) and the embedding $X\hookrightarrow L^r(M)$,
\[
|B_u(h,v)|
\leq
\|f_t(\cdot,u)\|_{L^s(M)}\|h\|_{L^r(M)}\|v\|_{L^r(M)}
\leq
C_S^2\|f_t(\cdot,u)\|_{L^s(M)}\|h\|_X\|v\|_X.
\]
Therefore, $B_u\in\mathcal L^2(X;\mathbb R)$.

For $h,v\in X$, Theorem~\ref{teo:b5-fundamental-calculo-riemann-banach}, applied in the second variable, gives
\[
f(x,u+h)-f(x,u)
=
\int_0^1f_t(x,u+\tau h)h\,d\tau
\]
for almost every $x\in M$. Thus,
\[
\bigl(DP(u+h)-DP(u)-B_u(h,\cdot)\bigr)(v)
=
\int_0^1\int_M
\bigl(f_t(x,u+\tau h)-f_t(x,u)\bigr)hv\,d\lambda_{\mathbf{g}}\,d\tau.
\]
Taking the supremum over $\|v\|_X\leq1$ and applying Hölder, we obtain
\[
\|DP(u+h)-DP(u)-B_u(h,\cdot)\|_{X'}
\leq
C_S^2\|h\|_X
\int_0^1
\|f_t(\cdot,u+\tau h)-f_t(\cdot,u)\|_{L^s(M)}\,d\tau.
\]
If $h_j\to0$ in $X$, then $u+\tau h_j\to u$ in $L^r(M)$ for each $\tau\in[0,1]$. Continuity of the Nemytski operator associated with $f_t$ implies that the integrand converges to zero. Moreover, the family $(u+\tau h_j)_{j,\tau}$ is bounded in $L^r(M)$, and the boundedness estimate of Corollary~\ref{cor:nemytski-crecimiento-variacional} provides a constant independent of $j$ and $\tau$ that dominates the integrand. By Theorem~\ref{convergencia dominada}, applied in $[0,1]$,
\[
\int_0^1
\|f_t(\cdot,u+\tau h_j)-f_t(\cdot,u)\|_{L^s(M)}\,d\tau
\longrightarrow0.
\]
Consequently,
\[
\lim_{h\to0}
\frac{\|DP(u+h)-DP(u)-B_u(h,\cdot)\|_{X'}}{\|h\|_X}=0.
\]
This proves that $DP$ is Fréchet differentiable and that $D^2P(u)=B_u$.

If $u_j\to u$ in $X$, then $u_j\to u$ in $L^r(M)$ and
\[
\|D^2P(u_j)-D^2P(u)\|_{\mathcal L^2(X;\mathbb R)}
\leq
C_S^2\|f_t(\cdot,u_j)-f_t(\cdot,u)\|_{L^s(M)}
\longrightarrow0.
\]
Therefore, $P\in C^2(X,\mathbb R)$. Since the second derivative of the quadratic term is the form $(v,w)\mapsto\displaystyle\int_M\langle\nabla v,\nabla w\rangle_{\mathbf{g}}\,d\lambda_{\mathbf{g}}$, the formula for $D^2I$ follows from $I=Q-P$.
\end{proof}

In the subcritical regime, the compact embedding of $W_0^{1,2}(M)$ into $L^r(M)$ makes the potential term stable under the weak convergences that occur in the direct method and in the Nehari method.

\chapter{The direct method in the calculus of variations}
\label{cap:metodo-directo-calculo-variaciones}

The direct method turns a minimization problem into a combination of three properties: coercivity, compactness, and lower semicontinuity. Coercivity keeps minimizing sequences in bounded sets; weak compactness allows a convergent subsequence to be extracted; and lower semicontinuity ensures that the limit does not increase the energy. In reflexive spaces, this scheme gives a general criterion for the existence of minimizers.

The distinction between weak compactness and weak convergence of sequences allows two versions of the method to be formulated. After establishing them, we will apply them to a sublinear elliptic problem. In that example, the Dirichlet condition controls the quadratic part of the energy, whereas compactness of a Sobolev embedding allows passage to the limit in the nonlinear term.

\begin{semblanzaHistorica}{The direct method and the power of weak convergence}
The direct method does not attempt to solve the Euler--Lagrange equation explicitly. It starts from a sequence whose energy approaches the infimum and seeks enough compactness to obtain a limit. Weak convergence is decisive because every bounded sequence in a reflexive Banach space admits a weakly convergent subsequence; lower semicontinuity prevents the energy from jumping upward at the limit. This combination, developed in the modern calculus of variations, turns structural properties of the space and the functional into an existence theorem.
\end{semblanzaHistorica}

\section{Lower semicontinuity}

A minimizing sequence often converges only weakly. For its limit to remain a minimizer, the energy must not increase on passing to the limit; this is precisely the role of lower semicontinuity.
\label{sec:semicontinuidad-inferior}

\begin{definition}\label{def:semicontinuidad-inferior}\index{lower semicontinuity}\glsadd{semicontinuidad-inferior}
 Let $X$ be a topological space. A function
$\varphi\colon X\longrightarrow\mathbb{R}$ is said to be lower semicontinuous if, for every $a\in\mathbb{R}$, the set $\{x\in X\mid\varphi(x)>a\}$ is open in $X$.
\end{definition}

\begin{example}\label{ej:continua-implica-semicontinua-inferior}
 Every continuous function $\varphi\colon X\longrightarrow\mathbb{R}$ is lower semicontinuous, since the preimage of any open set is an open set in $X$.
\end{example}

\begin{example}\label{ej:funcion-escalon-semicontinua-inferior}
 Consider the function $\varphi\colon \mathbb{R}\longrightarrow\mathbb{R}$ defined by

\[
\varphi(x) =
\begin{cases}
0, & \text{if } x \leq 0, \\
1, & \text{if } x > 0.
\end{cases}
\]

Then the function $\varphi$ is lower semicontinuous. Indeed, if $a<0$, the set $\{x\in\mathbb R\mid\varphi(x)>a\}$ is $\mathbb R$; if $0\leq a<1$, it is $(0,+\infty)$; and if $a\geq1$, it is empty. All three sets are open. The discontinuity at $0$ is compatible with lower semicontinuity because the value at that point agrees with the smaller of the two one-sided values.
\end{example}

\begin{lemma}\label{lem:criterio-secuencial-semicontinuidad-inferior}
Let $X$ be a first-countable topological space, and let $\psi\colon X\longrightarrow\mathbb{R}$ be a functional. Then $\psi$ is lower semicontinuous if and only if

\[
\psi(x) \leq \liminf_{n \to \infty} \psi(x_n),
\]

for every sequence $(x_n)_{n\in\mathbb N} \subset X$ such that $x_n \to x$ in $X$.
\end{lemma}
\begin{proof}
Suppose first that $\psi$ is lower semicontinuous. We wish to prove that $\displaystyle \psi(x) \leq \liminf_{n \to \infty} \psi(x_n)$ for every sequence $(x_n)_{n\in\mathbb N} \subset X$ such that $x_n \to x$. Let $\varepsilon > 0$. Consider the set
\[
V_\varepsilon := \{y \in X\mid \psi(y) > \psi(x) - \varepsilon\}
= \psi^{-1}\big((\psi(x)-\varepsilon, +\infty)\big).
\]
Since $\psi$ is lower semicontinuous, the set $V_\varepsilon$ is open. Since $x_n \to x$, there exists $n_0 \in \mathbb{N}$ such that $x_n \in V_\varepsilon$ for every $n \geq n_0$. Therefore,
\[
\psi(x_n) > \psi(x) - \varepsilon, \quad \text{for every } n \geq n_0,
\]
which implies $\displaystyle\liminf_{n \to \infty} \psi(x_n) \geq \psi(x) - \varepsilon$.

Since $\varepsilon > 0$ is arbitrary, it follows that $\displaystyle \psi(x) \leq \liminf_{n \to \infty} \psi(x_n)$.

Conversely, suppose that $\psi(x) \leq \displaystyle\liminf_{n \to \infty} \psi(x_n)$ for every sequence $(x_n)_{n\in\mathbb N} \subset X$ such that $x_n \to x$. We show that $\psi$ is lower semicontinuous. Let $a \in \mathbb{R}$ and consider the set $F := \{x \in X\mid \psi(x) \leq a\}$. We claim that $F$ is closed. Suppose, for a contradiction, that it is not. Then there exists $x_0 \in \overline{F} \setminus F$; that is, $x_0$ is a point of the closure of $F$ but $\psi(x_0) > a$. Since $X$ is first countable, $x_0$ has a countable base of open neighborhoods $(U_n)_{n\in\mathbb N}$. For each $n\in\mathbb N$, set $V_n:=\displaystyle\bigcap_{j=1}^{n}U_j$. Since $x_0\in\overline F$, we can choose $x_n\in F\cap V_n$. If $U$ is a neighborhood of $x_0$, there exists $j\in\mathbb N$ such that $U_j\subseteq U$; for every $n\geq j$, we have $x_n\in V_n\subseteq U_j\subseteq U$. Therefore, $x_n\to x_0$. By hypothesis, $\psi(x_0) \leq \displaystyle\liminf_{n \to \infty} \psi(x_n)$.

Since $x_n \in F$, we have $\psi(x_n) \leq a$ for every $n\in\mathbb N$, and hence $\displaystyle\liminf_{n \to \infty} \psi(x_n) \leq a$. Consequently, $\psi(x_0) \leq a$, contradicting $\psi(x_0) > a$. Therefore, $F$ is closed for every $a \in \mathbb{R}$, which implies that $\psi$ is lower semicontinuous.
\end{proof}
\begin{theorem}\label{teo:minimo-funcional-semicontinuo-compacto}\index{lower semicontinuity!existence of a minimum on compact spaces}
Let $X$ be a nonempty compact topological space, and let
$\varphi\colon X\longrightarrow\mathbb{R}$ be a lower semicontinuous functional. Then $\varphi$ is bounded below and attains its infimum; that is, there exists $u_0 \in X$ such that
\[
\varphi(u_0) = \inf_{x \in X} \varphi(x).
\]
\end{theorem}

\begin{proof}
Observe that $X = \displaystyle\bigcup_{n=1}^{\infty} \varphi^{-1}((-n,+\infty))$.

Since $\varphi$ is lower semicontinuous, each set
$\varphi^{-1}((-n,+\infty))$ is open in $X$.

Since $X$ is compact, the preceding cover admits a finite subcover. If $n_0\in\mathbb N$ is the largest of its indices, then $X = \displaystyle\bigcup_{n=1}^{n_0} \varphi^{-1}((-n,+\infty))$.

Therefore, for every $x \in X$ we have $\varphi(x) > -n_0$, which implies that $\varphi$ is bounded below.

Since $X\neq\varnothing$, the set $\varphi(X)$ is nonempty and bounded below. Completeness of $\mathbb R$ then provides a real number $c=\displaystyle\inf_{x\in X}\varphi(x)$.

Suppose, for a contradiction, that $\varphi(x) > c$ for every $x \in X$. Then $X = \displaystyle\bigcup_{n=1}^{\infty} \varphi^{-1}\left(\left(c + \displaystyle\frac{1}{n}, +\infty\right)\right)$.

Again, compactness provides a finite subcover. Taking $k\in\mathbb N$ to be the largest index of that subcover, we obtain $X = \displaystyle\bigcup_{n=1}^{k} \varphi^{-1}\left(\left(c + \displaystyle\frac{1}{n}, +\infty\right)\right)$.

It follows that $\varphi(x) > c + \displaystyle\frac{1}{k}$ for every $x\in X$, which implies that $c + \displaystyle\frac{1}{k}$ is a lower bound for $\varphi(X)$. This contradicts the definition of the infimum, since it would imply $c \geq c + \displaystyle\frac{1}{k}$. Therefore, there exists $u_0 \in X$ such that $\varphi(u_0) = \displaystyle\inf_{x \in X} \varphi(x)$.

\end{proof}

\begin{definition}\label{def:semicontinuidad-inferior-debil}\index{lower semicontinuity!weak}
Let $X$ be a normed space. A function $\varphi\colon X\longrightarrow\mathbb{R}$ is said to be weakly lower semicontinuous if it is lower semicontinuous with respect to the weak topology $\sigma(X,X')$. In particular, if a sequence $(x_n)_{n\in\mathbb N}\subseteq X$ satisfies $x_n\rightharpoonup x$ in $X$, then
\[
\varphi(x) \leq \liminf_{n \to \infty} \varphi(x_n).
\]
\end{definition}

The implication from lower semicontinuity to the inequality for sequences does not require first countability: the first part of the proof of Lemma~\ref{lem:criterio-secuencial-semicontinuidad-inferior} suffices. The converse does require additional justification when working with the weak topology.

\begin{definition}\label{def:semicontinuidad-inferior-debil-secuencial}
Let $X$ be a normed space. A functional $\varphi\colon X\longrightarrow\mathbb R$ is sequentially weakly lower semicontinuous if, for every sequence $(x_n)_{n\in\mathbb N}\subseteq X$ and every $x\in X$ such that $x_n\rightharpoonup x$, we have
\[
\varphi(x)\leq\liminf_{n\to\infty}\varphi(x_n).
\]
\end{definition}

The weak topology of an infinite-dimensional Banach space need not be first countable. We therefore keep the two notions distinct. In the applications, it will suffice to check the sequential condition and use a version of the direct method formulated for it.

\begin{theorem}[Direct method in the calculus of variations]\label{teo:metodo-directo-reflexivo}\index{direct method in the calculus of variations@direct method in the calculus of variations}\index{method!direct}\index{coercivity}\glsadd{metodo-directo}\glsadd{coercividad}
 Let $E$ be a reflexive Banach space, and let
$\varphi\colon E\longrightarrow\mathbb{R}$ be a functional satisfying:
\begin{itemize}
\item[(i)] $\varphi$ is weakly lower semicontinuous;
\item[(ii)] $\varphi$ is coercive, that is, $\varphi(u)\to+\infty$ as $\|u\|\to+\infty$: for every $a\in\mathbb R$ there exists $R_a>0$ such that $\varphi(u)>a$ for every $u\in E$ with $\|u\|\geq R_a$.
\end{itemize}

Then $\varphi$ is bounded below and attains its infimum; that is, there exists $u_0 \in E$ such that

\[
\varphi(u_0) = \inf_{u \in E} \varphi(u).
\]
\end{theorem}
\begin{proof}

Let $R > 0$ be such that
\[
\varphi(u) \geq \varphi(0), \quad \text{for every } u \in E \text{ with } \|u\| \geq R.
\]

Since $E$ is a reflexive Banach space, Theorem~\ref{teo: kakutani} implies that the closed ball $\overline{B}_{d_E}(0,R)=\{u\in E\mid\|u\|\leq R\}$, where $d_E$ is the metric induced by the norm, is weakly compact.

Moreover, the restriction of $\varphi$ to $\overline{B}_{d_E}(0,R)$ is weakly lower semicontinuous. Therefore, applying Theorem~\ref{teo:minimo-funcional-semicontinuo-compacto} shows that $\varphi$ is bounded below on $\overline{B}_{d_E}(0,R)$ and that there exists $u_0\in\overline{B}_{d_E}(0,R)$ such that $\varphi(u_0)=\displaystyle\inf_{u\in\overline{B}_{d_E}(0,R)}\varphi(u)$.

By the choice of $R$, for every $u \in E$ with $\|u\| \geq R$ we have $\varphi(u) \geq \varphi(0) \geq \varphi(u_0)$, which implies that $\displaystyle \varphi(u_0) = \inf_{u \in E} \varphi(u)$.

\end{proof}

Coercivity is equivalent to boundedness of every sublevel set $\{u\in E\mid\varphi(u)\leq a\}$, with $a\in\mathbb R$. Indeed, the radius $R_a$ in the definition encloses that sublevel set in a ball. Conversely, if the sublevel set is contained in a ball of radius $r_a$, any $R_a>r_a$ satisfies the coercivity condition. This observation allows us to work directly with minimizing sequences.

\begin{theorem}[Sequential form of the direct method]
\label{teo:metodo-directo-reflexivo-secuencial}
Let $E$ be a reflexive Banach space, and let $\varphi\colon E\longrightarrow\mathbb R$ be a coercive, sequentially weakly lower semicontinuous functional. Then $\varphi$ is bounded below and there exists $u_0\in E$ such that
\[
\varphi(u_0)=\inf_{u\in E}\varphi(u).
\]
\end{theorem}

\begin{proof}
First we rule out the possibility that the infimum is $-\infty$. If it were, for each $n\in\mathbb N$ we could choose $u_n\in E$ with $\varphi(u_n)\leq-n$. The entire sequence would lie in the sublevel set $\{u\in E\mid\varphi(u)\leq-1\}$, which is bounded by coercivity. Corollary~\ref{cor: sucesion acotada reflexivo tiene subsucesion debil} gives a subsequence $(u_{n_j})_{j\in\mathbb N}$ and a point $u\in E$ such that $u_{n_j}\rightharpoonup u$. Sequential semicontinuity would give
\[
\varphi(u)\leq\liminf_{j\to\infty}\varphi(u_{n_j})=-\infty,
\]
which is impossible because $\varphi$ takes real values. Consequently, $c:=\displaystyle\inf_{u\in E}\varphi(u)$ is finite; moreover, $c\leq\varphi(0)<+\infty$.

By the definition of the infimum, we can now choose, for each $n\in\mathbb N$, a point $v_n\in E$ such that
\[
c\leq\varphi(v_n)<c+\frac1n.
\]
This sequence lies in the sublevel set at height $c+1$ and is therefore bounded. Applying the cited corollary again, we obtain $u_0\in E$ and a subsequence $(v_{n_j})_{j\in\mathbb N}$ such that $v_{n_j}\rightharpoonup u_0$. Then
\[
c\leq\varphi(u_0)\leq\liminf_{j\to\infty}\varphi(v_{n_j})=c.
\]
The two inequalities force $\varphi(u_0)=c$.
\end{proof}

The extraction of weakly convergent subsequences rests on Eberlein--Šmulian and reflexivity, through the cited corollary in Appendix B. This is also the formulation of the direct method used in \cite[chap.~I, \S1, Theorem~1.2]{Struwe2008}.

Under the hypotheses of the last theorem, topological lower semicontinuity can also be recovered. Fix $a\in\mathbb R$ and set $K_a:=\{u\in E\mid\varphi(u)\leq a\}$. If $K_a$ is nonempty, every sequence in $K_a$ is bounded and admits a weakly convergent subsequence in $E$. Sequential semicontinuity ensures that its limit remains in $K_a$. Thus, $K_a$ is sequentially weakly compact and, by Theorem~\ref{teo: eberlein smulian}, weakly compact. Since the weak topology is Hausdorff, $K_a$ is weakly closed. The empty sublevel set is also closed, so all sublevel sets are closed. Proposition~\ref{prop: semicontinuidad inferior subniveles} then implies that $\varphi$ is weakly lower semicontinuous. This argument uses both coercivity and reflexivity; it has not assumed that the two notions of semicontinuity agree in an arbitrary weak topology.

We turn to an application in which the energy is not convex. Its quadratic part will be weakly lower semicontinuous, and a compact embedding will allow us to handle the remaining part.

\section{Solving a sublinear problem}
\label{sec:problema-sublineal-metodo-directo}

Let $(M,\mathbf{g})$ be a nonempty compact Riemannian manifold of dimension $n\geq1$, with smooth boundary, and suppose that every connected component of $M$ has nonempty boundary. We will work with real functions. The condition on the components allows the norm of a function with zero trace to be controlled by its gradient.

We will use the Poincaré inequality of Theorem~\ref{teo:poincare-W0-haces}, applied to the trivial real bundle with $p=2$. Its hypothesis on each component is precisely the one just imposed, and its range $1\leq p<\infty$ includes this exponent in every dimension $n\geq1$. The constant depends on $M$ and $\mathbf g$.

Consider the sublinear problem

\begin{equation}\label{eq:problema-sublineal-laplace}
\begin{cases}
-\Delta_{\mathbf{g}} u=|u|^{q-2}u & \text{in }M,\\
u=0 & \text{on }\partial M,
\end{cases}
\end{equation}
where $1<q<2$. We use the sign $\Delta_{\mathbf g}=\operatorname{div}_{\mathbf g}\nabla$, so that $-\Delta_{\mathbf g}$ is nonnegative under the Dirichlet condition. The expression $|t|^{q-2}t$ is defined to be zero when $t=0$; with this convention it is continuous on all of $\mathbb R$. The condition $u=0$ on $\partial M$ is understood in the trace sense.

The aim is to find $u\in W_0^{1,2}(M)$ such that

\begin{equation*}
 \int_M \langle \nabla u, \nabla v \rangle_{\mathbf{g}} \, d\lambda_{\mathbf{g}}
=
\int_M |u|^{q-2}u\, v \, d\lambda_{\mathbf{g}}
\end{equation*}
for every $v\in W_0^{1,2}(M)$. The associated energy functional is

\[
I(u) =
\frac{1}{2} \int_M |\nabla u|_{\mathbf{g}}^2 \, d\lambda_{\mathbf{g}}
-
\frac{1}{q} \int_M |u|^q \, d\lambda_{\mathbf{g}}.
\]

\begin{proposition}[Application to the sublinear problem]\label{prop:existencia-problema-sublineal}\index{application to the sublinear problem@application to the sublinear problem}
Suppose that every connected component of the compact manifold $M$ has nonempty smooth boundary. If $1<q<2$, the functional $I$ attains its global minimum on $W_0^{1,2}(M)$, and every minimizer is a nontrivial weak solution of problem \eqref{eq:problema-sublineal-laplace}.
\end{proposition}

\begin{proof}
The nonlinearity $f(x,t)=|t|^{q-2}t$ satisfies the hypotheses of Proposition~\ref{prop:funcional-energia-escalar-C1} with exponent $r=q$: it is continuous and its growth is $|t|^{q-1}$, with $1<q<2\leq2^*$. Hence, $I\in C^1(W_0^{1,2}(M),\mathbb R)$ and its critical points are precisely the weak solutions. More explicitly,
\[
DI(u)(v)=\int_M\langle\nabla u,\nabla v\rangle_{\mathbf g}\,d\lambda_{\mathbf g}-\int_M|u|^{q-2}u\,v\,d\lambda_{\mathbf g}
\]
for all $u,v\in W_0^{1,2}(M)$. The second term is finite by Hölder, since $|u|^{q-2}u\in L^{q/(q-1)}(M)$ and $v\in L^q(M)$.

Let $C_P$ be the Poincaré constant from Theorem~\ref{teo:poincare-W0-haces}. Since $M$ has finite volume, Hölder and Poincaré give
\[
\|u\|_{L^q(M)}\leq\lambda_{\mathbf g}(M)^{\frac1q-\frac12}\|u\|_{L^2(M)},
\qquad
\|u\|_{W^{1,2}(M)}^2\leq(1+C_P^2)\|\nabla u\|_{L^2(M,TM)}^2.
\]
Consequently, we can take $C_1:=\displaystyle\frac{1}{2(1+C_P^2)}$ and $C_2:=\displaystyle\frac1q\lambda_{\mathbf g}(M)^{1-q/2}$ to obtain
\[
I(u)\geq C_1\|u\|_{W^{1,2}(M)}^2-C_2\|u\|_{W^{1,2}(M)}^q.
\]
If $t:=\|u\|_{W^{1,2}(M)}>0$, the right-hand side of this estimate is $t^2(C_1-C_2t^{q-2})$. Since $q-2<0$, there exists $R>0$ such that $C_2t^{q-2}\leq C_1/2$ for every $t\geq R$. Thus, $I(u)\geq C_1\|u\|_{W^{1,2}(M)}^2/2$ as $\|u\|_{W^{1,2}(M)}\geq R$, proving coercivity.

To study passage to the limit, let $(u_j)_{j\in\mathbb N}\subseteq W_0^{1,2}(M)$ be a sequence such that $u_j\rightharpoonup u$ in that space. By Proposition~\ref{prop: convergencia debil acotada lsc}, there exists $R_0>0$ such that $\|u_j\|_{W^{1,2}(M)}\leq R_0$ for every $j\in\mathbb N$. The preceding estimate implies $I(u_j)\geq-C_2R_0^q$, whereas the definition of $I$ gives $I(u_j)\leq R_0^2/2$. Hence, $\ell:=\displaystyle\liminf_{j\to\infty}I(u_j)$ is a real number. Choose a subsequence $(u_{j_k})_{k\in\mathbb N}$ such that $I(u_{j_k})\to\ell$. Corollary~\ref{cor:rellich-descenso-un-orden-haces}, applied to the trivial real bundle with $m=1$ and $p=2$, gives compactness of \(W^{1,2}(M)\hookrightarrow L^2(M)\). Since
\(\lambda_{\mathbf g}(M)<\infty\), Hölder gives the continuous
embedding
\[
\|w\|_{L^q(M)}
\leq
\lambda_{\mathbf g}(M)^{\frac1q-\frac12}
\|w\|_{L^2(M)}
\qquad(1<q<2).
\]
By composition, \(W^{1,2}(M)\hookrightarrow L^q(M)\) is compact.
Passing to a subsequence of $(u_{j_k})_{k\in\mathbb N}$, which we continue to denote in the same way, we obtain strong convergence in $L^q(M)$. Its limit is $u$: the continuous embedding into $L^q(M)$ carries weak convergence of $u_{j_k}$ to weak convergence in $L^q(M)$, and every strong limit is also a weak limit. Uniqueness of the weak limit identifies the two limits. In particular, $\|u_{j_k}\|_{L^q(M)}\to\|u\|_{L^q(M)}$, and hence
\[
\int_M|u_{j_k}|^q\,d\lambda_{\mathbf{g}}\longrightarrow
\int_M|u|^q\,d\lambda_{\mathbf{g}}.
\]
The operator $u\mapsto\nabla u$ is linear and continuous from $W_0^{1,2}(M)$ to $L^2(M,TM)$, so $\nabla u_{j_k}\rightharpoonup\nabla u$. Weak lower semicontinuity of the norm, Proposition~\ref{prop: convergencia debil acotada lsc}, and convergence of the nonlinear term give
\[
I(u)\leq\liminf_{k\to\infty}\left(\frac12\|\nabla u_{j_k}\|_{L^2(M,TM)}^2-\frac1q\|u_{j_k}\|_{L^q(M)}^q\right)=\ell.
\]
Thus, $I$ is sequentially weakly lower semicontinuous. The space $W_0^{1,2}(M)$ is Hilbert, being a closed subspace of $W^{1,2}(M)$, and in particular is reflexive. Theorem~\ref{teo:metodo-directo-reflexivo-secuencial} then provides $u_0\in W_0^{1,2}(M)$ such that $I(u_0)=\displaystyle\inf_{u\in W_0^{1,2}(M)}I(u)$.

It remains to prove that $u_0\neq0$. Take $v\in C_c^\infty(\operatorname{Int}(M))\setminus\{0\}$. For $t>0$,
\[
I(tv)=\frac{t^2}{2}\int_M|\nabla v|_{\mathbf{g}}^2\,d\lambda_{\mathbf{g}}-\frac{t^q}{q}\int_M|v|^q\,d\lambda_{\mathbf{g}}.
\]
Write $A:=\displaystyle\int_M|\nabla v|_{\mathbf g}^2\,d\lambda_{\mathbf g}$ and $B:=\displaystyle\int_M|v|^q\,d\lambda_{\mathbf g}$. We have $B>0$, since $v\neq0$, and $A>0$ by Poincaré. The preceding equality becomes $I(tv)=t^q(At^{2-q}/2-B/q)$. Since $2-q>0$, any
\[
0<t_0<\left(\frac{2B}{qA}\right)^{\frac1{2-q}}
\]
satisfies $I(t_0v)<0=I(0)$. Consequently, $I(u_0)\leq I(t_0v)<0$ and $u_0\neq0$. For each $v\in W_0^{1,2}(M)$, the real map $t\mapsto I(u_0+tv)$ has a minimum at $t=0$, so its derivative there is zero. Thus, $DI(u_0)(v)=0$ for every $v\in W_0^{1,2}(M)$, which is the asserted weak solution identity.
\end{proof}

If only $\partial M\neq\varnothing$ were assumed, there could be a closed component $M_0$. For the function $u_t$ equal to $t$ on $M_0$ and zero on the other components, we would have
\[
I(u_t)=-\frac{\lambda_{\mathbf g}(M_0)}q|t|^q\longrightarrow-\infty
\qquad\text{as }|t|\to\infty.
\]
In that case the energy has no global minimum on $W_0^{1,2}(M)$. Nevertheless, a nontrivial weak solution can still be obtained if there is some component with boundary: apply the proposition on that component and extend the solution by zero to the remaining ones. The weak identity on $M$ reduces to the identity on the chosen component, since both the function and its gradient vanish on the others. This extension does not assert that the solution minimizes the energy on all of $M$.

\section{Exercises}
\label{sec:ejercicios-metodo-directo}

\begin{exercise}\label{ejer:minimizacion-conjunto-convexo-cerrado}
 Let $E$ be a reflexive Banach space, and let $C\subset E$ be a nonempty closed convex set.

Prove that if $\psi\colon E\longrightarrow\mathbb{R}$ is a weakly lower semicontinuous and coercive functional, then $\psi$ attains its infimum on $C$; that is, there exists $x_0 \in C$ such that

\[
\psi(x_0) = \inf_{x \in C} \psi(x).
\]

\emph{Hint.} Theorem~\ref{teo: cerradura convexa fuerte y debil} implies that $C$ is weakly closed. Repeat the sequential proof of the direct method with sequences contained in $C$; their weak limits will also belong to $C$. To bound the infimum above, use the value of $\psi$ at a fixed point of $C$, which exists because $C\neq\varnothing$.
\end{exercise}

\begin{exercise}\label{ejer:minimizacion-funcional-cuadratico-convexo}
 Let $E$ be a real Hilbert space, and let
$a\colon E\times E\longrightarrow\mathbb{R}$ be a continuous bilinear form satisfying the coercivity condition

\[
a(u,u) \geq \alpha \|u\|^2, \quad \text{for every } u \in E,
\]

for some $\alpha > 0$. Let also $f\colon E\longrightarrow\mathbb{R}$ be a continuous linear functional.

Given a nonempty closed convex set $C\subset E$, prove that there exists a unique $u_0 \in C$ such that

\[
\psi(u_0) = \inf_{u \in C} \psi(u),
\]

where the functional $\psi\colon E\longrightarrow\mathbb{R}$ is defined by

\[
\psi(u) = \frac{1}{2} a(u,u) - f(u).
\]

\emph{Hint.} It is not necessary to assume that $a$ is symmetric. Define $a_{\mathrm s}(u,v):=(a(u,v)+a(v,u))/2$ and observe that $a(u,u)=a_{\mathrm s}(u,u)$. The quadratic part is continuous and strictly convex; by Corollary~\ref{cor: funcional convexo semicontinuo es debilmente semicontinuo}, it is weakly lower semicontinuous. To make uniqueness precise, check that
\[
\psi((1-t)u+tv)=(1-t)\psi(u)+t\psi(v)-\frac{t(1-t)}2a(v-u,v-u)
\]
for all $u,v\in E$ and every $t\in(0,1)$. Also characterize the minimum by $a_{\mathrm s}(u_0,v-u_0)\geq f(v-u_0)$ for every $v\in C$.
\end{exercise}

\begin{exercise}\label{ejer:funcional-energia-laplaciano-cubico}
 Consider the problem

\begin{equation}\label{eq:ejercicio-laplaciano-cubico}
\begin{cases}
-\Delta u = u^2 & \text{in } \Omega, \\
u = 0 & \text{on } \partial \Omega,
\end{cases}
\end{equation}

where $\Omega\subset\mathbb R^3$ is a bounded domain. Work in $W_0^{1,2}(\Omega)$, and interpret the Dirichlet condition through the definition of this space as the closure of $C_c^\infty(\Omega)$ in $W^{1,2}(\Omega)$. This formulation does not require assigning a pointwise trace to an arbitrary boundary.

\begin{itemize}
\item[(a)] Define an energy functional associated with problem \eqref{eq:ejercicio-laplaciano-cubico} such that, if this functional is of class $C^1$, its critical points are weak solutions of \eqref{eq:ejercicio-laplaciano-cubico}.

\item[(b)] Prove that this functional is of class $C^1$.
\end{itemize}

\emph{Hint.} Use $F(t)=t^3/3$, whose derivative is $t^2$. The potential $|t|^3/3$ would have a different derivative. The embedding $W_0^{1,2}(\Omega)\hookrightarrow L^3(\Omega)$ allows the integral of the potential and its derivative to be controlled. Also check that this energy tends to $-\infty$ along $tv$ as $t\to+\infty$, when $v\in C_c^\infty(\Omega)\setminus\{0\}$ is nonnegative. This distinguishes the variational formulation from the hypotheses needed to obtain a global minimum.

\end{exercise}

\begin{exercise}\label{ejer:metodo-directo-problema-biharmonico}
Let $(M,\mathbf{g})$ be a nonempty compact Riemannian manifold of dimension $n\geq1$, with smooth boundary, and suppose that every connected component has nonempty boundary. Let $1<q<2$. Suppose that the elliptic estimate
\[
\|u\|_{H^2(M)}
\leq
C\|\Delta_{\mathbf{g}} u\|_{L^2(M)}
\]
holds on $H_0^2(M)$. Here $H_0^2(M)=W_0^{2,2}(M)$, and $\nu$ is the outward unit normal. The conditions $u=0$ and $\partial_\nu u=0$ are interpreted in the trace sense: Theorem~\ref{teo:caracterizacion-W0-trazas} requires vanishing of the traces up to order one, and tangential derivatives of the zero trace of $u$ also vanish. Consider the problem with clamped boundary conditions
\begin{equation}\label{eq:ejercicio-problema-biharmonico}
\begin{cases}
\Delta_{\mathbf{g}}^2u=|u|^{q-2}u&\text{in }M,\\
u=0,\quad\partial_\nu u=0&\text{on }\partial M.
\end{cases}
\end{equation}
\begin{enumerate}
\item Rigorously formulate the notion of a weak solution in $H_0^2(M)$.
\item Prove that the functional
\[
I(u)
:=
\frac{1}{2}\int_M|\Delta_{\mathbf{g}} u|^2\,d\lambda_{\mathbf{g}}
-
\frac{1}{q}\int_M|u|^q\,d\lambda_{\mathbf{g}}
\]
is well defined and belongs to $C^1(H_0^2(M),\mathbb R)$, and compute $DI(u)$.
\item Prove that $I$ is coercive and sequentially weakly lower semicontinuous, and deduce that it attains its global minimum.
\item Prove that the minimizer is nonzero and is a weak solution of \eqref{eq:ejercicio-problema-biharmonico}.
\end{enumerate}

\emph{Hint.} Use the elliptic estimate to work with the norm $\|\Delta_{\mathbf{g}} u\|_{L^2(M)}$, the compact embedding of $H_0^2(M)$ into $L^q(M)$, and the fact that $q<2$. To prove that the minimizer is nonzero, evaluate $I(tv)$ with $v\in C_c^\infty(\operatorname{Int}(M))\setminus\{0\}$ and small $t>0$.
\end{exercise}

\begin{exercise}\label{ejer:metodo-directo-problema-fraccionario}
Let $(M,\mathbf{g})$ be a nonempty compact Riemannian manifold of dimension $n\geq1$, with smooth boundary, such that every connected component has nonempty boundary. Let
$-\Delta_{\mathbf{g},D}$ be the nonnegative self-adjoint Dirichlet realization of the
Laplacian. Let $0<s<1$ and $1<q<2$, and define
$H_{0,D}^s(M):=\mathcal D((-\Delta_{\mathbf{g},D})^{s/2})$. Suppose it is a Hilbert
space under the norm
\[
\|u\|_{H_{0,D}^s(M)}
:=
\|(-\Delta_{\mathbf{g},D})^{\frac{s}{2}}u\|_{L^2(M)}
\]
and that the inclusion $H_{0,D}^s(M)\hookrightarrow L^q(M)$ is compact. Consider
\begin{equation}\label{eq:ejercicio-problema-fraccionario}
(-\Delta_{\mathbf{g},D})^su=|u|^{q-2}u\quad\text{in }M.
\end{equation}
Formulate the corresponding weak identity and prove, using the direct method, that there exists a nontrivial weak solution.

\emph{Hint.} Study the functional
\[
I(u)
:=
\frac{1}{2}\int_M|(-\Delta_{\mathbf{g},D})^{\frac{s}{2}}u|^2\,d\lambda_{\mathbf{g}}
-
\frac{1}{q}\int_M|u|^q\,d\lambda_{\mathbf{g}}.
\]
Prove separately that it is well defined, that it is of class $C^1$, that it is coercive, and that it is sequentially weakly lower semicontinuous. The assumed compactness gives strong convergence in $L^q(M)$ along a subsequence of any weakly convergent sequence in $H_{0,D}^s(M)$. Apply Theorem~\ref{teo:metodo-directo-reflexivo-secuencial}. The space $H_{0,D}^s(M)$ is defined through the spectral Dirichlet realization: no pointwise trace condition is added for every $s\in(0,1)$, nor is it identified, without another result, with the closure of $C_c^\infty(\operatorname{Int}(M))$ in a different fractional norm. In the weak formulation, the powers of order $s/2$ act on $u$ and on the test function; it is not initially required that $u$ belong to $\mathcal D((-\Delta_{\mathbf g,D})^s)$.
\end{exercise}

\begin{exercise}\label{ejer:metodo-directo-problema-kirchhoff}
Let $(M,\mathbf{g})$ be a nonempty compact Riemannian manifold of dimension $n\geq1$, with smooth boundary, such that every connected component has nonempty boundary. Let $1<q<2$, and let $\mathscr M\colon [0,+\infty)\longrightarrow(0,+\infty)$ be continuous and nondecreasing. Suppose that there exists $m_0>0$ such that $\mathscr M(t)\geq m_0$ for every $t\geq0$, and define
\[
\widehat{\mathscr M}(t)
:=
\int_0^t\mathscr M(s)\,ds.
\]
Consider the problem
\begin{equation}\label{eq:ejercicio-problema-kirchhoff}
\begin{cases}
-\mathscr M\!\left(\displaystyle\int_M|\nabla u|_{\mathbf{g}}^2\,d\lambda_{\mathbf{g}}\right)\Delta_{\mathbf{g}} u
=|u|^{q-2}u&\text{in }M,\\
u=0&\text{on }\partial M.
\end{cases}
\end{equation}
Prove that the functional
\[
I(u)
:=
\frac{1}{2}\widehat{\mathscr M}\!\left(\int_M|\nabla u|_{\mathbf{g}}^2\,d\lambda_{\mathbf{g}}\right)
-
\frac{1}{q}\int_M|u|^q\,d\lambda_{\mathbf{g}},
\qquad u\in W_0^{1,2}(M),
\]
is of class $C^1$, compute its derivative, and prove that it has a nonzero global minimizer, which is a weak solution of \eqref{eq:ejercicio-problema-kirchhoff}.

\emph{Hint.} Monotonicity of $\mathscr M$ implies that $\widehat{\mathscr M}$ is convex, whereas $\mathscr M\geq m_0$ gives $\widehat{\mathscr M}(t)\geq m_0t$. Moreover, $\widehat{\mathscr M}$ is increasing because $\mathscr M>0$. Combine this property with lower semicontinuity of $\|\nabla u\|_{L^2(M,TM)}^2$, and use the compact embedding $W_0^{1,2}(M)\hookrightarrow L^q(M)$ for the nonlinear term. Apply the chain rule to compute $DI$. To obtain a negative energy value, fix $v\neq0$ and observe that, by continuity of $\mathscr M$ at $0$, there exists $C_0>0$ such that $\widehat{\mathscr M}(t^2\|\nabla v\|_{L^2(M,TM)}^2)\leq C_0t^2\|\nabla v\|_{L^2(M,TM)}^2$ for every sufficiently small $t>0$. Then compare $t^2$ with $t^q$.
\end{exercise}

\begin{exercise}\label{ejer:metodo-directo-problema-integrodiferencial}
Let $(M,\mathbf{g})$ be a nonempty compact Riemannian manifold of positive dimension, with or without smooth boundary, and let $K\colon M\times M\longrightarrow(0,+\infty)$ be a measurable symmetric kernel. Suppose that
\[
\|u\|_{X_K}^2
:=
\int_M\int_M|u(x)-u(y)|^2K(x,y)\,d\lambda_{\mathbf{g}}(x)\,d\lambda_{\mathbf{g}}(y)
\]
is the square of a norm on a nonzero real Hilbert space $X_K(M)$ of measurable functions, identified when they agree almost everywhere, and that, for some $1<q<2$, the inclusion $X_K(M)\hookrightarrow L^q(M)$ is continuous, injective, and compact. The requirement that $\|\cdot\|_{X_K}$ be a norm excludes nonzero constants from $X_K(M)$, because their difference $u(x)-u(y)$ would be zero. We do not work with classes of functions modulo constants: the term $\displaystyle \int_M|u|^q\,d\lambda_{\mathbf g}$ would not be well defined on that quotient. Consider the bilinear form
\[
\mathcal B_K(u,v)
:=
\int_M\int_M
(u(x)-u(y))(v(x)-v(y))K(x,y)\,d\lambda_{\mathbf{g}}(x)\,d\lambda_{\mathbf{g}}(y).
\]
\begin{enumerate}
\item Define $\mathcal L_K\colon X_K(M)\longrightarrow X_K(M)'$ by $\langle\mathcal L_Ku,v\rangle:=\mathcal B_K(u,v)$ for all $u,v\in X_K(M)$, and formulate the weak equation $\mathcal L_Ku=|u|^{q-2}u$ as an equality in $X_K(M)'$. On the right-hand side, a function in $L^{q/(q-1)}(M)$ acts on $v\in X_K(M)$ by integration.
\item Prove that
\[
I(u)
:=
\frac{1}{2}\mathcal B_K(u,u)
-
\frac{1}{q}\int_M|u|^q\,d\lambda_{\mathbf{g}}
\]
belongs to $C^1(X_K(M),\mathbb R)$ and that its critical points are precisely the weak solutions.
\item Apply the direct method to prove the existence of a nontrivial weak solution.
\end{enumerate}

\emph{Hint.} Cauchy--Schwarz on the measure space $K(x,y)\,d\lambda_{\mathbf g}(x)\,d\lambda_{\mathbf g}(y)$ gives $|\mathcal B_K(u,v)|\leq\|u\|_{X_K}\|v\|_{X_K}$. The quadratic part is convex and weakly lower semicontinuous. For the nonlinear part, use the assumed compact embedding and differentiability of the functional $u\mapsto\|u\|_{L^q(M)}^q/q$. Nontriviality follows by comparing $I(tv)$ with $I(0)$ for a fixed element $v\in X_K(M)\setminus\{0\}$ and sufficiently small $t>0$. The equation is asserted only against test functions in $X_K(M)$; if this space incorporates a constraint, an identity against all smooth test functions does not follow without further argument.
\end{exercise}

\begin{exercise}\label{ejer:metodo-directo-problema-p-laplaciano}
Let $(M,\mathbf{g})$ be a nonempty compact Riemannian manifold of dimension $n\geq1$ with smooth boundary, and suppose that every connected component has nonempty boundary. Let $1<p<\infty$ and $1<q<p$, and consider
\[
\begin{cases}
-\Delta_{p,\mathbf{g}} u=|u|^{q-2}u&\text{in }M,\\
u=0&\text{on }\partial M,
\end{cases}
\qquad
\Delta_{p,\mathbf{g}} u
:=
\operatorname{div}_{\mathbf{g}}\bigl(|\nabla u|_{\mathbf{g}}^{p-2}\nabla u\bigr).
\]
Formulate the weak problem in $W_0^{1,p}(M)$ and prove that it has a nontrivial weak solution by minimizing
\[
I(u)
:=
\frac{1}{p}\int_M|\nabla u|_{\mathbf{g}}^p\,d\lambda_{\mathbf{g}}
-
\frac{1}{q}\int_M|u|^q\,d\lambda_{\mathbf{g}}.
\]

\emph{Hint.} Corollary~\ref{cor:potencia-norma-Lp-C1-haces} and the chain rule give differentiability of both terms. Use Theorem~\ref{teo:poincare-W0-haces} for the trivial real bundle. Corollary~\ref{cor:rellich-descenso-un-orden-haces}, with $m=1$, gives compactness $W^{1,p}(M)\hookrightarrow L^p(M)$ for every $1\leq p<\infty$. Since $q<p$ and the volume is finite, its composition with $L^p(M)\hookrightarrow L^q(M)$ gives the required compact embedding. Convexity of the principal term provides its weak lower semicontinuity; $1<p<\infty$ ensures reflexivity of $W_0^{1,p}(M)$. Apply Theorem~\ref{teo:metodo-directo-reflexivo-secuencial} and compare $t^p$ with $t^q$ to prove nontriviality. In the expression for the $p$-Laplacian, define $|\xi|^{p-2}\xi=0$ when $\xi=0$, and all derivatives in the weak problem are understood in the sense of distributions.
\end{exercise}

\chapter{Palais--Smale, pseudogradients, and deformations}
\label{cap:palais-smale-deformacion}

The direct method produces critical points by minimization. Minimax
methods require a different mechanism: the sublevel sets of a functional must be
deformed within bands containing no critical points. On a Hilbert
space this deformation can be constructed, under sufficient regularity, using
the negative gradient flow. On a general Banach space there is no
canonical identification with the dual, so the gradient is replaced by
a descent field retaining only the indispensable estimates.

We will construct these fields and prove the deformation lemma first on
Banach spaces and then on complete Banach--Finsler manifolds. The
Palais--Smale condition will appear exactly where it is necessary
to prevent an almost critical sequence from escaping without converging. The exposition follows Struwe~\cite[chap.~II, \S3, pp.~81--87]{Struwe2008}: the linear case corresponds to Lemma~II.3.2 and Theorem~II.3.4, pp.~81--84; the intrinsic version to \S\S II.3.7--II.3.10 and Theorem~II.3.11, pp.~85--87.

\begin{semblanzaHistorica}{Palais and Smale: compactness in the search for saddle points}
Finite-dimensional Morse theory uses gradient flow to deform sublevel sets without crossing critical points. Palais and Smale identified a compactness condition that preserves this strategy in infinite-dimensional spaces. A sequence with controlled energy and small derivative must have a convergent subsequence. When this holds, topological changes in the sublevel sets are tied to actual critical points, and minimax arguments acquire analytic content~\cite{PalaisSmale1964,Palais1970}.
\end{semblanzaHistorica}

\section{Sublevel sets, critical points, and the Palais--Smale condition}
\label{sec:palais-smale-banach}

Let $X$ be a real Banach space, and let $I\in C^1(X,\mathbb R)$. We denote
by $X'$ the continuous dual of $X$ and by $\|DI(u)\|_{X'}$ the dual norm of the
derivative of $I$ at $u$.

\begin{definition}[Critical sets and sublevel sets]
\label{def:conjuntos-criticos-subniveles}
\index{critical point}
\index{critical set}
\index{sublevel set}
The critical set of $I$, the critical set at level $\beta\in\mathbb R$,
and the open sublevel set at height $a\in\mathbb R$ are, respectively,
\[
 \operatorname{Crit}(I):=\{u\in X\mid DI(u)=0\},
 \qquad
 K_\beta:=\{u\in X\mid I(u)=\beta,\ DI(u)=0\},
\]
and
\[
 I^{<a}:=\{u\in X\mid I(u)<a\}.
\]
The closed sublevel set is denoted by
$I^{\leq a}:=\{u\in X\mid I(u)\leq a\}$.
\end{definition}

\begin{definition}[Palais--Smale condition]
\label{def:palais-smale-banach}\glsadd{palais-smale}
\index{Palais--Smale condition}
A sequence $(u_j)\subset X$ is a \emph{Palais--Smale sequence} for
$I$ if $(I(u_j))$ is bounded and
\[
 \|DI(u_j)\|_{X'}\longrightarrow0.
\]
The functional $I$ satisfies the Palais--Smale condition, abbreviated
\emph{(PS)}, if every Palais--Smale sequence has a norm-convergent
subsequence.

For $\beta\in\mathbb R$, we say that $I$ satisfies \emph{$(PS)_\beta$} if
every sequence $(u_j)$ such that
\[
 I(u_j)\longrightarrow\beta,
 \qquad
 \|DI(u_j)\|_{X'}\longrightarrow0
\]
has a convergent subsequence in $X$.
\end{definition}

The global condition (PS) implies $(PS)_\beta$ for every $\beta$. The levelwise
version is exactly the compactness involved in a deformation
around a fixed height. These formulations and their relation to the
original condition of Palais and Smale are discussed in
\cite[chap.~II, \S2, pp.~77--80]{Struwe2008}.

For $\delta,\rho>0$ set
\begin{align*}
 N_{\beta,\delta}
 &:=
 \left\{u\in X\middle|
 |I(u)-\beta|<\delta,\
 \|DI(u)\|_{X'}<\delta
 \right\},\\
 U_{\beta,\rho}
 &:=
 \left\{u\in X\middle|
 \operatorname{dist}(u,K_\beta)<\rho
 \right\}.
\end{align*}
If $K_\beta=\varnothing$, adopt the convention
$U_{\beta,\rho}=\varnothing$.

\begin{lemma}[Local consequence of Palais--Smale]
\label{lem:consecuencia-ps-nivel-banach}
Suppose that $I$ satisfies $(PS)_\beta$. Then:
\begin{enumerate}[label=(\alph*)]
\item $K_\beta$ is compact;
\item $(U_{\beta,\rho})_{\rho>0}$ is a neighborhood base of
$K_\beta$;
\item $(N_{\beta,\delta})_{\delta>0}$ is a neighborhood base of
$K_\beta$.
\end{enumerate}
In particular, if $N$ is an open neighborhood of $K_\beta$, there exist
$\delta,\rho\in(0,1]$ such that
\begin{equation}
 N\supseteq U_{\beta,2\rho}
 \supseteq U_{\beta,\rho}
 \supseteq N_{\beta,\delta}.
\label{eq:cadena-vecindades-ps-banach}
\end{equation}
\end{lemma}

\begin{proof}
Every sequence in $K_\beta$ is a $(PS)_\beta$ sequence. By hypothesis it has
a convergent subsequence, and continuity of $I$ and $DI$ shows that its
limit belongs to $K_\beta$. Since $X$ is metrizable, $K_\beta$ is compact.

Let $N$ be an open neighborhood of $K_\beta$. If no $U_{\beta,\rho}$
were contained in $N$, there would exist $\rho_j\to 0^{+}$ and
$u_j\in U_{\beta,\rho_j}\setminus N$. For each $j$ we could choose
$v_j\in K_\beta$ with $\|u_j-v_j\|<2\rho_j$. Compactness of $K_\beta$
would allow us to assume that $v_j\to v\in K_\beta$; then also $u_j\to v$,
contradicting openness of $N$.

If no $N_{\beta,\delta}$ were contained in $N$, we could choose
$u_j\in N_{\beta,\frac{1}{j}}\setminus N$. This would be a $(PS)_\beta$ sequence.
A subsequence would converge to a point of $K_\beta\subset N$, a contradiction.
This proves the three assertions. Finally, choose first
$\rho\in(0,1]$ so that $U_{\beta,2\rho}\subseteq N$ and then
$\delta\in(0,1]$ so that
$N_{\beta,\delta}\subseteq U_{\beta,\rho}$.
\end{proof}

\begin{corollary}[Uniform separation of the derivative from zero]
\label{cor:separacion-derivada-ps-banach}
If the hypotheses of the preceding lemma hold and $N$ is an open neighborhood
of $K_\beta$, there exists $\delta>0$ such that
\[
 |I(u)-\beta|<\delta,
 \quad u\notin N
 \quad\Longrightarrow\quad
 \|DI(u)\|_{X'}\geq\delta.
\]
\end{corollary}

\begin{proof}
Take $\delta$ so that $N_{\beta,\delta}\subseteq N$. If the first two
conditions hold and $\|DI(u)\|_{X'}<\delta$, then
$u\in N_{\beta,\delta}\subseteq N$, a contradiction.
\end{proof}

\section{Construction of pseudogradients on Banach spaces}
\label{sec:pseudogradientes-banach}

Write
\[
 \widetilde X:=X\setminus\operatorname{Crit}(I).
\]
This set is open because $DI$ is continuous.

\begin{definition}[Pseudogradient vector field]
\label{def:pseudogradiente-banach}\glsadd{pseudogradiente}
\index{pseudogradient}
A field $v\colon \widetilde X\longrightarrow X$ is a \emph{pseudogradient} for $I$ if it is
locally Lipschitz and, for every $u\in\widetilde X$, satisfies
\begin{align}
 \|v(u)\|
 &<2\min\{\|DI(u)\|_{X'},1\},
 \label{eq:pseudogradiente-cota-banach}\\
 DI(u)[v(u)]
 &>\min\{\|DI(u)\|_{X'},1\}\,\|DI(u)\|_{X'}.
 \label{eq:pseudogradiente-descenso-banach}
\end{align}
\end{definition}

The second inequality expresses that $-v$ is a descent direction. The
first simultaneously provides control near critical points and
a global speed bound. Replacing the gradient by a controlled descent
direction is the functional-analytic mechanism used in
\cite{Palais1970}.

\begin{theorem}[Existence of pseudogradients]
\label{teo:existencia-pseudogradiente-banach}
Every functional $I\in C^1(X,\mathbb R)$ admits a pseudogradient vector field on
$\widetilde X$.
\end{theorem}

\begin{proof}
Fix $u\in\widetilde X$ and write
\[
 q_u:=\|DI(u)\|_{X'}>0,
 \qquad
 m_u:=\min\{q_u,1\}.
\]
By the definition of the dual norm, there exists $z_u\in X$, with $\|z_u\|=1$, such
that, after changing its sign if necessary,
\[
 DI(u)[z_u]>\frac34q_u.
\]
Define $w_u:=\frac32m_u z_u$. Then
\[
 \|w_u\|=\frac32m_u<2m_u,
 \qquad
 DI(u)[w_u]>\frac98m_uq_u>m_uq_u.
\]
Since $DI$ is continuous and the inequalities are strict, there exists an
open neighborhood $W_u\subseteq\widetilde X$ of $u$ such that, for every
$y\in W_u$,
\begin{align}
 \|w_u\|&<2\min\{\|DI(y)\|_{X'},1\},
 \label{eq:direccion-local-cota}\\
 DI(y)[w_u]&>
 \min\{\|DI(y)\|_{X'},1\}\,\|DI(y)\|_{X'}.
 \label{eq:direccion-local-descenso}
\end{align}

The metric space $\widetilde X$ is paracompact by Theorem~\ref{teo:paracompacidad-espacios-metricos}. Choose a locally finite open
refinement $(A_i)_{i\in J}$ of the cover
$(W_u)_{u\in\widetilde X}$ and, for each $i$, a point $u_i$ such that
$A_i\subseteq W_{u_i}$. Applying the locally finite shrinking in
Lemma~\ref{lem:encogimiento-localmente-finito-banach} gives an open
cover $(B_i)_{i\in J}$ of $\widetilde X$ such that
\[
 \overline{B_i}^{\,\widetilde X}\subseteq A_i
 \qquad(i\in J).
\]
This shrinking is essential: the distance functions used below
have support contained in the closure of $B_i$, not
necessarily in $B_i$.

For $y\in\widetilde X$ set
\[
 r_i(y):=\operatorname{dist}
 \bigl(y,\widetilde X\setminus B_i\bigr).
\]
If $B_i=\widetilde X$, interpret $r_i\equiv1$. Each $r_i$ is
Lipschitz, is positive precisely on $B_i$, and
\[
 \operatorname{supp}(r_i)
 \subseteq\overline{B_i}^{\,\widetilde X}
 \subseteq A_i.
\]
The sum $\displaystyle r:=\displaystyle\sum_{i\in J}r_i$ is locally finite and strictly positive.
Therefore,
\[
 \varphi_i:=\frac{r_i}{r}
\]
defines a locally Lipschitz partition of unity on $\widetilde X$, with
$\operatorname{supp}(\varphi_i)\subseteq A_i$.

Define
\[
 v(y):=\sum_{i\in J}\varphi_i(y)w_{u_i}.
\]
The sum is locally finite, so $v$ is locally Lipschitz. If
$\varphi_i(y)>0$, then $y\in A_i\subseteq W_{u_i}$ and the
inequalities \eqref{eq:direccion-local-cota} and
\eqref{eq:direccion-local-descenso} hold for $w_{u_i}$ at $y$. Taking the
convex combination gives
\eqref{eq:pseudogradiente-cota-banach} and
\eqref{eq:pseudogradiente-descenso-banach}. This proves the theorem.
\end{proof}

\section{The deformation lemma on a Banach space}
\label{sec:lema-deformacion-banach}

Two formulations appear in the literature and should be distinguished. The
first starts from a uniform derivative bound throughout an energy
band and requires no compactness hypothesis. The second, due to
Struwe, uses Palais--Smale to obtain this bound outside a neighborhood of the
critical set. We will state and prove both, then make their logical
relationship precise. The first form is customary in quantitative
treatments of minimax; see, for example, \cite{Willem1996}.

\subsection{A quantitative lemma without Palais--Smale}

\begin{theorem}[Quantitative deformation lemma]
\label{teo:lema-deformacion-cuantitativo-banach}
Let $X$ be a real Banach space, let $I\in C^1(X,\mathbb R)$, and fix
$\beta\in\mathbb R$, $a>0$, and $\sigma>0$. Suppose that
\begin{equation}
 |I(u)-\beta|\leq 2a
 \quad\Longrightarrow\quad
 \|DI(u)\|_{X'}\geq\sigma.
\label{eq:banda-regular-deformacion-cuantitativa}
\end{equation}
Then, for each $\varepsilon\in(0,a)$, there exists a continuous map
\[
 H\colon[0,1]\times X\longrightarrow X,
 \qquad H_t(u):=H(t,u),
\]
satisfying:
\begin{enumerate}[label=(\roman*)]
\item $H_0=\operatorname{id}_X$ and each $H_t$ is a homeomorphism of $X$;
\item if $|I(u)-\beta|\geq2a$, then $H_t(u)=u$ for every $t\in[0,1]$;
\item $t\mapsto I(H_t(u))$ is nonincreasing for every $u\in X$;
\item
\[
 H_1\bigl(I^{\leq\beta+\varepsilon}\bigr)
 \subseteq I^{\leq\beta-\varepsilon};
\]
\item if $\mu:=\min\{\sigma,1\}\sigma$, then
\[
 \|H_t(u)-u\|\leq\frac{4\varepsilon}{\mu}
 \qquad(t\in[0,1],\ u\in X).
\]
\end{enumerate}
In particular, the conclusion is quantitative both in the decrease of
energy and in the displacement of points.
\end{theorem}

\begin{proof}
Let $v$ be the pseudogradient from
Theorem~\ref{teo:existencia-pseudogradiente-banach}. Choose a Lipschitz
function $\chi\colon\mathbb R\to[0,1]$ such that
\[
 \chi(s)=1\quad\text{if }|s-\beta|\leq a,
 \qquad
 \chi(s)=0\quad\text{if }|s-\beta|\geq2a.
\]
For example, one can interpolate linearly on each of the two intermediate
intervals. Define
\begin{equation}
 e(u):=
 \begin{cases}
 -\chi(I(u))v(u),&u\notin\operatorname{Crit}(I),\\
 0,&u\in\operatorname{Crit}(I).
 \end{cases}
\label{eq:campo-lema-cuantitativo}
\end{equation}

We carefully verify that $e$ is locally Lipschitz on all of $X$. This property
is immediate outside $\operatorname{Crit}(I)$, since $I$ is locally
Lipschitz, $\chi$ is Lipschitz, and $v$ is locally Lipschitz. By
\eqref{eq:banda-regular-deformacion-cuantitativa} there are no critical points in the
closed band $|I-\beta|\leq2a$. Hence, if $u$ is critical, then
$|I(u)-\beta|>2a$; continuity of $I$ shows that $\chi\circ I$ vanishes
in a neighborhood of $u$, and there $e\equiv0$. Moreover, by
\eqref{eq:pseudogradiente-cota-banach},
\begin{equation}
 \|e(u)\|\leq2
 \qquad(u\in X).
\label{eq:velocidad-lema-cuantitativo}
\end{equation}

Theorem~\ref{teo:flujo-global-campo-acotado-banach} provides a global
flow $\Phi\colon X\times\mathbb R\to X$ for $e$. Uniqueness of
solutions implies that each $\Phi_s$ is a homeomorphism with inverse
$\Phi_{-s}$, and also that a point where $e$ vanishes remains fixed. By the
chain rule,
\begin{equation}
 \frac{d}{ds}I(\Phi_s(u))
 =-\chi(I(\Phi_s(u)))
 DI(\Phi_s(u))[v(\Phi_s(u))]
 \leq0.
\label{eq:descenso-lema-cuantitativo}
\end{equation}
At critical points, the right-hand side of \eqref{eq:descenso-lema-cuantitativo} is interpreted as zero. In particular, energy does not increase.

Set
\[
 \mu:=\min\{\sigma,1\}\sigma>0,
 \qquad
 T:=\frac{2\varepsilon}{\mu},
 \qquad
 H_t(u):=\Phi_{tT}(u).
\]
Properties (i)--(iii) follow from the preceding argument. The speed bound
\eqref{eq:velocidad-lema-cuantitativo} gives
\[
 \|H_t(u)-u\|
 \leq\int_0^{tT}\|e(\Phi_s(u))\|\,ds
 \leq2tT
 \leq\frac{4\varepsilon}{\mu},
\]
which is (v).

Only (iv) remains to be proved. Let $u\in I^{\leq\beta+\varepsilon}$. If
$I(u)\leq\beta-\varepsilon$, monotonicity suffices. Suppose,
otherwise, that $I(u)>\beta-\varepsilon$ and
$I(\Phi_T(u))>\beta-\varepsilon$. Then, for every $s\in[0,T]$,
\[
 \beta-\varepsilon
 < I(\Phi_s(u))
 \leq I(u)
 \leq\beta+\varepsilon.
\]
Since $\varepsilon<a$, along the entire trajectory we have
$\chi(I(\Phi_s(u)))=1$ and, by
\eqref{eq:banda-regular-deformacion-cuantitativa},
$\|DI(\Phi_s(u))\|_{X'}\geq\sigma$. The descent property of the
pseudogradient then implies
\[
 \frac{d}{ds}I(\Phi_s(u))\leq-\mu.
\]
Integrating between $0$ and $T$ gives
\[
 I(\Phi_T(u))
 \leq I(u)-\mu T
 \leq\beta+\varepsilon-2\varepsilon
 =\beta-\varepsilon,
\]
a contradiction. Consequently,
$H_1(u)=\Phi_T(u)\in I^{\leq\beta-\varepsilon}$.
\end{proof}

Hypothesis \eqref{eq:banda-regular-deformacion-cuantitativa} says that the entire
band is regular, with the derivative uniformly bounded away from zero. It does not
assert compactness of sequences with bounded energy. This is why
the preceding lemma applies even when (PS) fails, provided that
the slope bound can be verified by another argument.

\subsection{Struwe's general lemma}

The following is the quantitative form of Theorem~II.3.4 in
\cite[pp.~83--84]{Struwe2008}. The proof uses continuation of the flow
in both time directions and estimates how long a trajectory
remains in the region of uniform descent, with the constants in the
statement.

\begin{theorem}[Struwe's deformation lemma]
\label{teo:lema-deformacion-struwe-banach}
Let $X$ be a real Banach space, and let $I\in C^1(X,\mathbb R)$ be a functional
satisfying (PS). Fix $\beta\in\mathbb R$, $\varepsilon_0>0$, and an
open neighborhood $N$ of $K_\beta$. Then there exist
$\varepsilon\in(0,\varepsilon_0)$ and a continuous flow
\[
 \Phi\colon X\times\mathbb R\longrightarrow X,
 \qquad
 \Phi_t(u):=\Phi(u,t),
\]
with the following properties:
\begin{enumerate}[label=(\roman*)]
\item $\Phi_0=\operatorname{id}_X$ and each $\Phi_t$ is a homeomorphism, with
$\Phi_t^{-1}=\Phi_{-t}$;
\item $\Phi_{t+s}=\Phi_t\circ\Phi_s$ for every $s,t\in\mathbb R$;
\item if $u\in\operatorname{Crit}(I)$ or
$|I(u)-\beta|\geq\varepsilon_0$, then $\Phi_t(u)=u$ for every
$t\in\mathbb R$;
\item $t\mapsto I(\Phi_t(u))$ is nonincreasing for every $u\in X$;
\item
\[
 \Phi_1\bigl(I^{<\beta+\varepsilon}\setminus N\bigr)
 \subseteq I^{<\beta-\varepsilon},
 \qquad
 \Phi_1\bigl(I^{<\beta+\varepsilon}\bigr)
 \subseteq I^{<\beta-\varepsilon}\cup N.
\]
\end{enumerate}
In particular, the restriction of $(\Phi_t)_{t\in\mathbb R}$ to $t\geq0$ is
the sublevel deformation semiflow.
\end{theorem}

\begin{proof}
By Lemma~\ref{lem:consecuencia-ps-nivel-banach} we can choose
$\delta,\rho\in(0,1]$ such that
\begin{equation}
 N\supseteq U_{\beta,2\rho}
 \supseteq U_{\beta,\rho}
 \supseteq N_{\beta,\delta}.
\label{eq:cadena-deformacion-banach}
\end{equation}
Let $\eta\colon X\longrightarrow[0,1]$ be a locally Lipschitz function such that
\begin{equation}
 \eta=0\quad\text{in }N_{\beta,\frac{\delta}{2}},
 \qquad
 \eta=1\quad\text{in }X\setminus N_{\beta,\delta}.
\label{eq:corte-eta-deformacion-banach}
\end{equation}
It can be constructed without any additional regularity of the functional.
Indeed,
\[
 C_{\beta,\frac{\delta}{2}}:=
 \left\{u\in X\middle|
 |I(u)-\beta|\leq\frac\delta2,\quad
 \|DI(u)\|_{X'}\leq\frac\delta2
 \right\},
 \qquad
 A:=C_{\beta,\frac{\delta}{2}},
 \qquad
 B:=X\setminus N_{\beta,\delta}
\]
are closed by continuity of $I$ and
$u\mapsto\|DI(u)\|_{X'}$. Moreover,
\[
 A\cap B=\varnothing,
 \qquad
 N_{\beta,\frac{\delta}{2}}\subseteq A\subseteq N_{\beta,\delta}.
\]
When $A$ and $B$ are nonempty, define
\[
 \eta(u):=
 \frac{d_A(u)}{q(u)},
 \qquad
 d_A(u):=\operatorname{dist}(u,A),
 \quad
 d_B(u):=\operatorname{dist}(u,B),
 \quad
 q(u):=d_A(u)+d_B(u).
\]
The functions $d_A$ and $d_B$ are $1$-Lipschitz. Since $A$ and $B$ are closed
and disjoint, $q(u)>0$ for every $u$: equality $q(u)=0$ would imply
$u\in A\cap B$. Fixing $u_0\in X$ and writing $q_0:=q(u_0)>0$, if
$\|u-u_0\|<\frac{q_0}{4}$, then
\[
 q(u)\geq q_0-|d_A(u)-d_A(u_0)|-|d_B(u)-d_B(u_0)|
 >\frac{q_0}{2}.
\]
Therefore, the denominator is locally bounded away from zero. For $u$ and
$\widetilde u$ in that ball, write
$a=d_A(u)$, $\widetilde a=d_A(\widetilde u)$,
$q=q(u)$, and $\widetilde q=q(\widetilde u)$. Then
\[
 \left|\frac a q-\frac{\widetilde a}{\widetilde q}\right|
 \leq\frac{|a-\widetilde a|}{q}
 +\frac{\widetilde a\,|q-\widetilde q|}{q\widetilde q}.
\]
Here $\widetilde a$ is locally bounded, whereas $d_A$, $d_B$, and
$q$ are Lipschitz. Thus, $\eta$ is locally Lipschitz. On $A$ we have
$d_A=0$, and on $B$ we have $d_B=0$; disjointness and closedness show,
respectively, that the other summand is positive. Consequently,
$\eta=0$ on $A$ and $\eta=1$ on $B$, proving
\eqref{eq:corte-eta-deformacion-banach}.

If $A=\varnothing\neq B$, take $\eta\equiv1$. If
$B=\varnothing$, take $\eta\equiv0$, including the case where both are
empty. Since $N_{\beta,\frac{\delta}{2}}\subseteq A$, these conventions
preserve the two required values in
\eqref{eq:corte-eta-deformacion-banach}.

Set
\[
 a:=\min\{\varepsilon_0,\frac{\delta}{4}\}
\]
and define the Lipschitz cutoff $\chi\colon \mathbb R\longrightarrow[0,1]$ by
\[
 \chi(s):=
 \begin{cases}
 1,
 & |s-\beta|\leq \frac{a}{2},\\[2mm]
 2-\dfrac{2|s-\beta|}{a},
 & \frac{a}{2}<|s-\beta|<a,\\[2mm]
 0,
 & |s-\beta|\geq a.
 \end{cases}
\]
Let $v$ be a pseudogradient provided by
Theorem~\ref{teo:existencia-pseudogradiente-banach}. Define the truncated
field
\begin{equation}
 e(u):=
 \begin{cases}
 -\eta(u)\chi(I(u))v(u),
 &u\notin\operatorname{Crit}(I),\\
 0,
 &u\in\operatorname{Crit}(I).
 \end{cases}
\label{eq:campo-truncado-deformacion-banach}
\end{equation}

The field $e$ is locally Lipschitz. On $\widetilde X$, the
pseudogradient $v$ is so by construction, $\eta$ is locally Lipschitz, and
$\chi\circ I$ is as well, since $\chi$ is Lipschitz and $I$ is of class
$C^1$ and therefore locally Lipschitz. On a ball where all three
factors are bounded and have finite Lipschitz constants, the estimate
\[
 \|abc-a'b'c'\|
 \leq |a-a'|\,|b|\,\|c\|
      +|a'|\,|b-b'|\,\|c\|
      +|a'|\,|b'|\,\|c-c'\|
\]
for $a,a',b,b'\in\mathbb R$ and $c,c'\in X$, applied to
$a=\eta(u)$, $b=\chi(I(u))$, and $c=v(u)$, proves that their product is
locally Lipschitz. If $u$ is critical and
$|I(u)-\beta|<\frac{\delta}{2}$, the cutoff $\eta$ vanishes in a neighborhood of $u$.
If $|I(u)-\beta|\geq\frac{\delta}{2}$, then $|I(u)-\beta|>a$ and
$\chi\circ I$ vanishes in a neighborhood of $u$. Consequently, $e$ is
identically zero near every critical point.
Moreover,
\begin{equation}
 \|e(u)\|\leq2
 \qquad(u\in X).
\label{eq:cota-campo-truncado-banach}
\end{equation}

Applying Theorem~\ref{teo:flujo-global-campo-acotado-banach} to the autonomous
problem
\begin{equation}
 \frac{d}{dt}\Phi_t(u)=e(\Phi_t(u)),
 \qquad
 \Phi_0(u)=u,
\label{eq:cauchy-deformacion-banach}
\end{equation}
gives a solution defined for all $t\in\mathbb R$, continuous with respect
to $(u,t)$. Uniqueness gives
\[
 \Phi_{t+s}=\Phi_t\circ\Phi_s,
 \qquad
 \Phi_t\circ\Phi_{-t}=\Phi_0=\operatorname{id}_X.
\]
Thus $\Phi_t^{-1}=\Phi_{-t}$ and each $\Phi_t$ is a homeomorphism.

If $u$ is critical or $|I(u)-\beta|\geq\varepsilon_0$, then $e(u)=0$;
uniqueness implies that the solution starting at $u$ is constant. By the chain
rule and the pseudogradient properties,
\begin{align}
 \frac{d}{dt}I(\Phi_t(u))
 &=DI(\Phi_t(u))
 \left[\frac{d}{dt}\Phi_t(u)\right]
 \notag\\
 &=-\eta(\Phi_t(u))\chi(I(\Phi_t(u)))
 DI(\Phi_t(u))[v(\Phi_t(u))]
 \leq0.
\label{eq:identidad-descenso-deformacion-banach}
\end{align}
At a critical point the right-hand side is interpreted as zero. This proves
the first four properties.

Now choose
\begin{equation}
 0<\varepsilon<
 \min\left\{\frac a2,\frac{\rho\delta^2}{4}\right\}.
\label{eq:eleccion-epsilon-deformacion-banach}
\end{equation}
Let $u\in I^{<\beta+\varepsilon}$ and suppose that
$I(\Phi_1(u))\geq\beta-\varepsilon$. Monotonicity implies, for
$0\leq t\leq1$,
\[
 \beta-\varepsilon
 \leq I(\Phi_t(u))
 <\beta+\varepsilon.
\]
Hence $\chi(I(\Phi_t(u)))=1$. Define
\[
 A_u:=\{t\in[0,1]\mid
 \Phi_t(u)\notin N_{\beta,\delta}\}.
\]
If $t\in A_u$, then
$|I(\Phi_t(u))-\beta|<\varepsilon<\delta$, so exclusion from
$N_{\beta,\delta}$ forces
$\|DI(\Phi_t(u))\|_{X'}\geq\delta$. Since $\delta\leq1$ and $\eta=1$ outside
$N_{\beta,\delta}$, from
\eqref{eq:identidad-descenso-deformacion-banach} we obtain
\[
 \frac{d}{dt}I(\Phi_t(u))\leq-\delta^2
 \qquad(t\in A_u).
\]
Integrating,
\begin{equation}
 I(\Phi_1(u))
 <\beta+\varepsilon-\delta^2|A_u|,
\label{eq:descenso-medida-tiempo-banach}
\end{equation}
where $|A_u|$ denotes Lebesgue measure.

Suppose that $u\notin N$ or $\Phi_1(u)\notin N$.
If $K_\beta=\varnothing$, choose $\delta$ with
$N_{\beta,\delta}=\varnothing$; then $A_u=[0,1]$.
If $K_\beta\neq\varnothing$, set
$h(t):=\operatorname{dist}(\Phi_t(u),K_\beta)$. The field bound gives
\[
 \|\Phi_t(u)-\Phi_s(u)\|
 \leq2|t-s|,
 \qquad |h(t)-h(s)|\leq2|t-s|.
\]
By \eqref{eq:cadena-deformacion-banach}, $u\notin N$ implies
$h(0)\geq2\rho$. Therefore,
\[
 h(t)\geq2\rho-2t\geq\rho
 \quad(0\leq t\leq\rho/2),
 \qquad [0,\rho/2]\subseteq A_u.
\]
If $\Phi_1(u)\notin N$, then $h(1)\geq2\rho$ and, similarly,
\[
 h(t)\geq2\rho-2(1-t)\geq\rho
 \quad(1-\rho/2\leq t\leq1),
 \qquad [1-\rho/2,1]\subseteq A_u.
\]
These inclusions use $N_{\beta,\delta}\subseteq U_{\beta,\rho}$.
Since $\rho\leq1$, in every case we obtain
\begin{equation}
 |A_u|\geq\frac\rho2.
\label{eq:tiempo-fuera-Nbd-banach}
\end{equation}

Estimates \eqref{eq:descenso-medida-tiempo-banach},
\eqref{eq:tiempo-fuera-Nbd-banach}, and
\eqref{eq:eleccion-epsilon-deformacion-banach} yield
\[
 I(\Phi_1(u))
 <\beta+\varepsilon-\frac{\rho\delta^2}{2}
 <\beta-\varepsilon,
\]
contradicting the preceding assumption. If
$u\in I^{<\beta+\varepsilon}\setminus N$, apply the argument to the initial
endpoint and obtain the first inclusion in (v). For an arbitrary
$u\in I^{<\beta+\varepsilon}$, either $\Phi_1(u)\in N$ or the
same argument applied to the final endpoint gives
$\Phi_1(u)\in I^{<\beta-\varepsilon}$. This proves the second inclusion.
\end{proof}

\begin{remark}[Levelwise hypotheses]
\label{obs:deformacion-ps-nivel-banach}
Condition (PS) was used only in
Lemma~\ref{lem:consecuencia-ps-nivel-banach}. Therefore, in the preceding
theorem it suffices to assume $(PS)_\beta$. If $K_\beta=\varnothing$, there exists
$\delta>0$ such that $N_{\beta,\delta}=\varnothing$; in this case we may take
$N=\varnothing$ and obtain a uniform decrease in energy in a
band around $\beta$.
\end{remark}

\subsection{Relationship between the two deformation lemmas}

The difference between the two statements lies in the source of the lower bound
for $\|DI\|$. The quantitative lemma assumes it; Struwe's lemma
deduces it from $(PS)_\beta$ outside a neighborhood of $K_\beta$. The next
result formalizes this relationship and specifies exactly which part of
Palais--Smale is used.

\begin{proposition}[From compactness to a slope bounded away from zero]
\label{prop:relacion-lemas-deformacion}
Let $I\in C^1(X,\mathbb R)$ and let $\beta\in\mathbb R$.
\begin{enumerate}[label=(\alph*)]
\item If $I$ satisfies $(PS)_\beta$ and $K_\beta=\varnothing$, there exist
$a,\sigma>0$ such that
\[
 |I(u)-\beta|\leq2a
 \quad\Longrightarrow\quad
 \|DI(u)\|_{X'}\geq\sigma.
\]
Therefore, Theorem~\ref{teo:lema-deformacion-cuantitativo-banach} uniformly deforms
the sublevel set just above $\beta$ into one
just below it.
\item If $I$ satisfies $(PS)_\beta$ and $N$ is an open neighborhood of
$K_\beta$, there exist $a,\sigma>0$ such that
\[
 |I(u)-\beta|\leq2a,\qquad u\notin N
 \quad\Longrightarrow\quad
 \|DI(u)\|_{X'}\geq\sigma.
\]
Struwe's lemma is the relative version of the same mechanism: an additional
cutoff stops the flow near $K_\beta$ and preserves uniform
decrease outside $N$.
\end{enumerate}
\end{proposition}

\begin{proof}
For (a), suppose there are no such $a$ and $\sigma$. For each
$j\in\mathbb N$ we could choose $u_j\in X$ with
\[
 |I(u_j)-\beta|\leq\frac{2}{j},
 \qquad
 \|DI(u_j)\|_{X'}<\frac1j.
\]
The sequence $(u_j)$ would be a $(PS)_\beta$ sequence. A subsequence
would converge to some $u\in X$, and continuity of $I$ and $DI$ would give
$I(u)=\beta$ and $DI(u)=0$, contradicting $K_\beta=\varnothing$.

For (b), Corollary~\ref{cor:separacion-derivada-ps-banach} provides
$\delta>0$ such that
\[
 |I(u)-\beta|<\delta,\qquad u\notin N
 \quad\Longrightarrow\quad
 \|DI(u)\|_{X'}\geq\delta.
\]
It suffices to take, for example, $a=\delta/4$ and $\sigma=\delta$. The proof of
Theorem~\ref{teo:lema-deformacion-struwe-banach} then shows why it does not
suffice to turn off the field inside $N$: the nested neighborhoods
$U_{\beta,\rho}\subset U_{\beta,2\rho}$ are used to control also the
trajectories entering or leaving $N$.
\end{proof}

Thus, the statements are not equivalent without additional hypotheses. If
$K_\beta\neq\varnothing$, the global hypothesis of the quantitative lemma is
impossible, because $DI$ vanishes on $K_\beta$. If $K_\beta=\varnothing$,
however, $(PS)_\beta$ turns the quantitative lemma into the case of Struwe's lemma
needed for minimax contradiction arguments.

\section{The intrinsic Banach--Finsler formulation}
\label{sec:palais-smale-finsler}

Now consider a reversible Banach--Finsler manifold without boundary $M$ with the
regularity and local compatibility of norms from
Section~\ref{sec:banach-finsler} in
Chapter~\ref{cap:geometria-hilbert-finsler}. Denote by $\|\xi\|_p$ the
norm of $\xi\in T_pM$, by $d$ the length distance, and by
$\|\alpha\|_{p,*}$ the dual norm of
$\alpha\in T'_pM$ defined in
Definition~\ref{def:norma-dual-finsler}. In particular,
\[
 \|\alpha\|_{p,*}
 =\sup\left\{|\alpha(\xi)|\middle|
 \xi\in T_pM,\ \|\xi\|_p\leq1\right\}.
\]
Reversibility is the property of the structure $F$ specified in
Definition~\ref{def:estructura-banach-finsler}. When class
$C^{1,1}$ is required, we will use the atlas of
Definition~\ref{def:variedad-banach-C11}, in accordance with
Definition~\ref{def:finsler-C11}. Locally Lipschitz fields and their
flows are understood as in
Remark~\ref{obs:campos-flujos-atlas-C11}.
Since all arguments apply componentwise, we assume below
that $M$ is connected.

\begin{definition}[Intrinsic Palais--Smale condition]
\label{def:palais-smale-finsler}
\index{Palais--Smale condition!on a Finsler manifold}
Let $M$ be a reversible Banach--Finsler manifold without boundary, and let
$I\in C^1(M,\mathbb R)$. A sequence $(p_j)\subset M$ is a
Palais--Smale sequence if $(I(p_j))$ is bounded and
\[
 \|dI_{p_j}\|_{p_j,*}\longrightarrow0.
\]
The functional satisfies (PS) if every Palais--Smale sequence has a
convergent subsequence in $M$. It satisfies $(PS)_\beta$ if the same conclusion
holds for every sequence such that
\[
 I(p_j)\longrightarrow\beta,
 \qquad
 \|dI_{p_j}\|_{p_j,*}\longrightarrow0.
\]
\end{definition}

Define intrinsically
\begin{align*}
 K_\beta
 &:=
 \{p\in M\mid I(p)=\beta,\ dI_p=0\},\\
 N_{\beta,\delta}
 &:=
 \{p\in M\mid |I(p)-\beta|<\delta,
 \ \|dI_p\|_{p,*}<\delta\},\\
 U_{\beta,\rho}
 &:=
 \{p\in M\mid d(p,K_\beta)<\rho\}.
\end{align*}
We also retain the notation
$I^{<a}:=\{p\in M\mid I(p)<a\}$ for open sublevel sets and
$I^{\leq a}:=\{p\in M\mid I(p)\leq a\}$ for closed ones.
Local compatibility of the Finsler structure implies that $d$ induces the
topology of $M$ and that $p\mapsto\|dI_p\|_{p,*}$ is continuous. Hence the
argument of Lemma~\ref{lem:consecuencia-ps-nivel-banach}, which uses only
metrizability, continuity, and the definition of $(PS)_\beta$, gives the following
result.

\begin{lemma}[Neighborhoods of the critical set on a Finsler manifold]
\label{lem:consecuencia-ps-nivel-finsler}
Let $M$ be a reversible Banach--Finsler manifold without boundary, and let
$I\in C^1(M,\mathbb R)$. If $I$ satisfies $(PS)_\beta$, then $K_\beta$ is compact and the families
$(U_{\beta,\rho})_{\rho>0}$ and
$(N_{\beta,\delta})_{\delta>0}$ are neighborhood bases of
$K_\beta$. In particular, for every open neighborhood $N$ of $K_\beta$,
there exist $\delta,\rho\in(0,1]$ such that
\[
 N\supseteq U_{\beta,2\rho}
 \supseteq U_{\beta,\rho}
 \supseteq N_{\beta,\delta}.
\]
\end{lemma}

\begin{proof}
A sequence in $K_\beta$ is a $(PS)_\beta$ sequence, and by continuity
every subsequential limit remains in $K_\beta$; therefore,
$K_\beta$ is compact. If a neighborhood $N$ contained no
$N_{\beta,\delta}$, a choice
$p_j\in N_{\beta,\frac{1}{j}}\setminus N$ would give a $(PS)_\beta$ sequence whose
convergent subsequence would have a limit in $K_\beta\subset N$, a contradiction.
The assertion for $U_{\beta,\rho}$ follows from compactness of $K_\beta$ and
the same metric argument: for points at decreasing distance from
$K_\beta$, choose nearby points in $K_\beta$ and extract a convergent
subsequence. The final chain is obtained by choosing first $\rho$ and then
$\delta$.
\end{proof}

\begin{definition}[Finsler pseudogradient]
\label{def:pseudogradiente-finsler}
Let $M$ be a reversible Banach--Finsler manifold without boundary, and let
$I\in C^1(M,\mathbb R)$. Let
\[
 \widetilde M:=\{p\in M\mid dI_p\neq0\}.
\]
A locally Lipschitz vector field
$\mathbf{v}\colon \widetilde M\longrightarrow TM$ is a
\emph{Finsler pseudogradient} of $I$ if $\mathbf{v}(p)\in T_pM$ and
\begin{align*}
 \|\mathbf{v}(p)\|_p
 &<2\min\{\|dI_p\|_{p,*},1\},\\
 dI_p[\mathbf{v}(p)]
 &>\min\{\|dI_p\|_{p,*},1\}\,\|dI_p\|_{p,*}
\end{align*}
for every $p\in\widetilde M$.
\end{definition}

\begin{theorem}[Existence of a Finsler pseudogradient]
\label{teo:existencia-pseudogradiente-finsler}
Let $M$ be a paracompact Banach--Finsler manifold of class $C^{1,1}$.
Then every $I\in C^1(M,\mathbb R)$ admits a Finsler pseudogradient on
$\widetilde M$.
\end{theorem}

\begin{proof}
Fix $p\in\widetilde M$ and write
$q_p:=\|dI_p\|_{p,*}$ and $m_p:=\min\{q_p,1\}$. By the definition of the dual
norm, there exists $\xi_p\in T_pM$, $\|\xi_p\|_p=1$, such that
\[
 dI_p[\xi_p]>\frac34q_p.
\]
Therefore,
\[
 \left\|\frac32m_p\xi_p\right\|_p
 =\frac32m_p<2m_p,
 \qquad
 dI_p\left[\frac32m_p\xi_p\right]
 >\frac98m_pq_p>m_pq_p.
\]
Choose a chart $(U,\varphi)$ of the $C^{1,1}$ atlas about $p$,
with $U\subseteq\widetilde M$, and the tangent trivialization
$\tau_q:=d\varphi_q\colon T_qM\longrightarrow X$. Define
\[
 \mathbf{w}_p(q)
 :=\tau_q^{-1}\tau_p\left(\frac32m_p\xi_p\right),
 \qquad q\in U.
\]
Its components in the chart $\varphi$ are constant, and
$\mathbf{w}_p(p)=\frac32m_p\xi_p$. By admissibility of the Finsler
structure, Definition~\ref{def:estructura-banach-finsler}, and comparison
\eqref{eq:comparacion-local-finsler-C}, after shrinking $U$ there exist
constants $c_p,C_p>0$ such that
\[
 c_p\|\tau_q(z)\|_X
 \leq \|z\|_q
 \leq C_p\|\tau_q(z)\|_X
\]
for every point $q$ in the domain and every $z\in T_qM$. The functions
$q\mapsto\|\mathbf{w}_p(q)\|_q$ and
$q\mapsto dI_q[\mathbf{w}_p(q)]$ are continuous by continuity of
$F$, $dI$, and the local field in this trivialization. Moreover,
$q\mapsto\|dI_q\|_{q,*}$ is continuous by
Proposition~\ref{prop:dual-finsler-coordenadas}, as is its minimum
with $1$. Since both inequalities at $p$ are strict,
we can shrink the domain again to a neighborhood $W_p$ on which
$\mathbf{w}_p$ satisfies
\begin{align*}
 \|\mathbf{w}_p(q)\|_q
 &<2\min\{\|dI_q\|_{q,*},1\},\\
 dI_q[\mathbf{w}_p(q)]
 &>\min\{\|dI_q\|_{q,*},1\}\,\|dI_q\|_{q,*}.
\end{align*}
In the chart $\varphi$, the field has constant components and is
locally Lipschitz. Remark~\ref{obs:campos-flujos-atlas-C11}
shows that this property is preserved in all charts of the
$C^{1,1}$ atlas.

The open set $\widetilde M$, with the restricted Finsler distance, is metrizable by Theorem~\ref{teo:distancia-banach-finsler-topologia} and paracompact by Theorem~\ref{teo:paracompacidad-espacios-metricos}. Choose a locally finite open refinement $(A_i)_{i\in J}$ of $(W_p)_{p\in\widetilde M}$ and, by Lemma~\ref{lem:encogimiento-localmente-finito-banach}, a shrinking $(B_i)_{i\in J}$ with $\overline{B_i}^{\,\widetilde M}\subseteq A_i$. Define
\[
 r_i(q):=d\bigl(q,\widetilde M\setminus B_i\bigr),
 \qquad
 r(q):=\sum_{i\in J}r_i(q),
 \qquad
 \varphi_i(q):=\frac{r_i(q)}{r(q)}.
\]
If $B_i=\widetilde M$, take $r_i\equiv1$. Each $r_i$ is Lipschitz with respect to $d$ and positive precisely on $B_i$. The sum is locally finite and positive because $(B_i)_{i\in J}$ covers $\widetilde M$. Near each point, $r$ is bounded away from zero; hence the quotients $\varphi_i$ are locally $d$-Lipschitz. Theorem~\ref{teo:distancia-banach-finsler-topologia} transfers this regularity to the charts. Moreover,
\[
 \sum_{i\in J}\varphi_i=1,
 \qquad
 \operatorname{supp}(\varphi_i)
 \subseteq\overline{B_i}^{\,\widetilde M}\subseteq A_i.
\]
For each $i$ choose $p_i$ with $A_i\subseteq W_{p_i}$ and set
\[
 \mathbf{v}(q):=\sum_{i\in J}\varphi_i(q)\mathbf{w}_{p_i}(q).
\]
In this formula, each product $\varphi_i\mathbf{w}_{p_i}$ is extended
by zero outside $W_{p_i}$. This extension is locally Lipschitz on
$\widetilde M$: inside $W_{p_i}$ it is a product of locally Lipschitz functions and fields,
which are locally bounded in coordinates; outside
$W_{p_i}$ it is identically zero on a neighborhood, since
$\operatorname{supp}(\varphi_i)\subseteq A_i\subseteq W_{p_i}$ and the
support is closed in $\widetilde M$.
Each sum is finite on a neighborhood of the point under consideration. Therefore,
$\mathbf{v}$ is
locally Lipschitz; since it is a convex combination of directions satisfying
the two inequalities on the corresponding support, it is a
Finsler pseudogradient.
\end{proof}

Completeness enters at a single point: it turns uniform speed control
into continuation of integral curves.

\begin{proposition}[Bounded fields on complete Finsler manifolds]
\label{prop:flujo-global-campo-acotado-finsler}
Let $(M,d)$ be a reversible complete Banach--Finsler manifold without boundary of class
$C^{1,1}$, and let $\mathbf{e}$ be a vector field locally Lipschitz in charts. If
there exists $C>0$ such that
\[
 \|\mathbf{e}(p)\|_p\leq C
 \qquad(p\in M),
\]
then all its maximal integral curves are defined on
$\mathbb R$. The flow $\Phi\colon M\times\mathbb R\longrightarrow M$ is continuous, satisfies
\[
 \Phi_{t+s}=\Phi_t\circ\Phi_s,
 \qquad
 \Phi_t^{-1}=\Phi_{-t},
\]
and each $\Phi_t$ is a homeomorphism.
\end{proposition}

\begin{proof}
Remark~\ref{obs:campos-flujos-atlas-C11} provides local existence,
uniqueness, and continuous dependence in the chosen atlas.
Let $\gamma\colon (a,b)\longrightarrow M$ be a maximal integral curve. For
$a<s<t<b$, the definition of the length distance gives
\begin{equation}
 d(\gamma(s),\gamma(t))
 \leq L(\gamma\restriction_{[s,t]})
 =\int_s^t\|\mathbf{e}(\gamma(\tau))\|_{\gamma(\tau)}\,d\tau
 \leq C|t-s|.
\label{eq:cauchy-flujo-finsler-completo}
\end{equation}
If $b<+\infty$, the estimate shows that $\gamma(t)$ is a
Cauchy family as $t\to b^{-}$. Completeness provides
$p_b\in M$ such that $\gamma(t)\to p_b$. Since the Finsler distance induces the
manifold topology, we can work in a chart about $p_b$.
Local existence for differential equations in Banach spaces produces an
integral curve starting at $p_b$ at time $b$. Continuity of $\mathbf{e}$ and
the integral equation allow it to be concatenated with $\gamma$, and local uniqueness
shows that this concatenation extends the original curve. This contradicts
maximality. Thus, $b=+\infty$. Applying the same argument as
$t\to a^{+}$ gives $a=-\infty$.

Continuous dependence on the initial value produces a continuous flow.
Uniqueness of integral curves gives the group law and, taking $s=-t$,
the identity $\Phi_t^{-1}=\Phi_{-t}$. Hence each $\Phi_t$ is a
homeomorphism. Observe that the argument uses metric completeness at
\eqref{eq:cauchy-flujo-finsler-completo}, not compactness of a closed
ball.
\end{proof}

\begin{theorem}[Banach--Finsler deformation lemma]
\label{teo:lema-deformacion-struwe-finsler}
Let $M$ be a reversible, complete, paracompact Banach--Finsler manifold without boundary, of
class $C^{1,1}$, satisfying the hypotheses of
Theorem~\ref{teo:existencia-pseudogradiente-finsler}. Let
$I\in C^1(M,\mathbb R)$ be a functional satisfying (PS). Given
$\beta\in\mathbb R$, $\varepsilon_0>0$, and an open neighborhood $N$ of
$K_\beta$, there exist $\varepsilon\in(0,\varepsilon_0)$ and a continuous flow
$\Phi\colon M\times\mathbb R\longrightarrow M$ such that:
\begin{enumerate}[label=(\roman*)]
\item each $\Phi_t$ is a homeomorphism,
$\Phi_t^{-1}=\Phi_{-t}$, and
$\Phi_{t+s}=\Phi_t\circ\Phi_s$;
\item $\Phi_t(p)=p$ if $dI_p=0$ or
$|I(p)-\beta|\geq\varepsilon_0$;
\item $I(\Phi_t(p))$ is nonincreasing in $t$;
\item
\[
 \Phi_1\bigl(I^{<\beta+\varepsilon}\setminus N\bigr)
 \subseteq I^{<\beta-\varepsilon},
 \qquad
 \Phi_1\bigl(I^{<\beta+\varepsilon}\bigr)
 \subseteq I^{<\beta-\varepsilon}\cup N.
\]
\end{enumerate}
Here too it suffices to assume $(PS)_\beta$.
\end{theorem}

\begin{proof}
Lemma~\ref{lem:consecuencia-ps-nivel-finsler} provides
$\delta,\rho\in(0,1]$ such that
\[
 N\supseteq U_{\beta,2\rho}
 \supseteq U_{\beta,\rho}
 \supseteq N_{\beta,\delta}.
\]
Take
\[
 a:=\min\{\varepsilon_0,\frac{\delta}{4}\}
\]
and construct the spatial cutoff entirely from the Finsler distance. Let
\[
 C^M_{\beta,\frac{\delta}{2}}:=
 \left\{p\in M\middle|
 |I(p)-\beta|\leq\frac\delta2,\quad
 \|dI_p\|_{p,*}\leq\frac\delta2
 \right\},
 \qquad
 A_M:=C^M_{\beta,\frac{\delta}{2}},
 \qquad
 B_M:=M\setminus N_{\beta,\delta}.
\]
Continuity of $I$ and $p\mapsto\|dI_p\|_{p,*}$ shows that $A_M$ and
$B_M$ are closed. Moreover,
\[
 A_M\cap B_M=\varnothing,
 \qquad
 N_{\beta,\frac{\delta}{2}}\subseteq A_M\subseteq N_{\beta,\delta}.
\]
If both sets are nonempty, set
\begin{equation}
 \eta(p):=\frac{d_A^M(p)}{q_M(p)},
 \quad
 d_A^M(p):=d(p,A_M),
 \quad
 d_B^M(p):=d(p,B_M),
 \quad
 q_M(p):=d_A^M(p)+d_B^M(p).
\label{eq:corte-eta-deformacion-finsler}
\end{equation}
The functions $d_A^M$ and $d_B^M$ are $1$-Lipschitz with respect to $d$.
Closedness and disjointness imply $q_M(p)>0$ for every $p$. Fixing
$p_0\in M$ and writing $q_0:=q_M(p_0)>0$, if $d(p,p_0)<\frac{q_0}{4}$, then
\[
 q_M(p)\geq q_0
 -|d_A^M(p)-d_A^M(p_0)|-|d_B^M(p)-d_B^M(p_0)|
 >\frac{q_0}{2}.
\]
Thus the denominator is locally bounded away from zero. For $p$ and
$\widetilde p$ in that ball, with
$a=d_A^M(p)$, $\widetilde a=d_A^M(\widetilde p)$,
$q=q_M(p)$, and $\widetilde q=q_M(\widetilde p)$, we have
\[
 \left|\frac a q-\frac{\widetilde a}{\widetilde q}\right|
 \leq\frac{|a-\widetilde a|}{q}
 +\frac{\widetilde a\,|q-\widetilde q|}{q\widetilde q}.
\]
The factor $\widetilde a$ is locally bounded, so $\eta$ is
locally $d$-Lipschitz. Local equivalence between $d$ and the norm of a
Finsler chart gives the required local regularity in charts. Since
$d_A^M=0$ on $A_M$ and $d_B^M=0$ on $B_M$, while the other summand is
positive, we obtain $\eta=0$ on $A_M$ and $\eta=1$ on $B_M$.

If $A_M=\varnothing\neq B_M$, take $\eta\equiv1$. If
$B_M=\varnothing$, take $\eta\equiv0$, including the case where both are
empty. In every case, $\eta$ equals zero on
$N_{\beta,\frac{\delta}{2}}$ and one outside $N_{\beta,\delta}$.

Let $\chi\colon\mathbb R\longrightarrow[0,1]$ be the Lipschitz function already
defined in the Banach proof: it equals one on
$[\beta-\frac{a}{2},\beta+\frac{a}{2}]$ and zero outside $(\beta-a,\beta+a)$.

Let $\mathbf{v}$ be the pseudogradient from
Theorem~\ref{teo:existencia-pseudogradiente-finsler} and define
\[
 \mathbf{e}(p):=
 \begin{cases}
 -\eta(p)\chi(I(p))\mathbf{v}(p),&dI_p\neq0,\\
 0,&dI_p=0.
 \end{cases}
\]
On the open set $\widetilde M$, the three factors $\eta$, $\chi\circ I$, and
$\mathbf{v}$ are locally Lipschitz in charts and locally bounded; the product
inequality used in the Banach proof shows that $\mathbf{e}$ is locally
Lipschitz there. Critical points remain to be checked. If $dI_p=0$ and
$|I(p)-\beta|<\frac{\delta}{2}$, then $p$ belongs to the open set
$N_{\beta,\frac{\delta}{2}}$, where $\eta=0$. If $dI_p=0$ and
$|I(p)-\beta|\geq\frac{\delta}{2}$, then
$|I(p)-\beta|\geq\frac{\delta}{2}>a$; by continuity, $\chi\circ I$ vanishes on
a neighborhood of $p$. In both cases $\mathbf{e}$ is identically zero on a
neighborhood of the critical point. Thus $\mathbf e$ is locally
Lipschitz on all of $M$. Moreover,
\begin{equation}
 \|\mathbf{e}(p)\|_p\leq2.
\label{eq:cota-campo-truncado-finsler}
\end{equation}
Proposition~\ref{prop:flujo-global-campo-acotado-finsler} gives a global flow
of homeomorphisms defined in both time directions. In particular, completeness is used
explicitly in the estimate
\begin{equation}
 d(\Phi_s(p),\Phi_t(p))
 \leq\int_s^t\|\mathbf{e}(\Phi_\tau(p))\|_{\Phi_\tau(p)}\,d\tau
 \leq2|t-s|.
\label{eq:velocidad-flujo-deformacion-finsler}
\end{equation}

The intrinsic chain rule gives
\begin{align*}
 \frac{d}{dt}I(\Phi_t(p))
 &=dI_{\Phi_t(p)}\bigl[\partial_t\Phi_t(p)\bigr]\\
 &=-\eta(\Phi_t(p))\chi(I(\Phi_t(p)))
 dI_{\Phi_t(p)}[\mathbf{v}(\Phi_t(p))]
 \leq0.
\end{align*}
Choose
\[
 0<\varepsilon<
 \min\left\{\frac a2,\frac{\rho\delta^2}{4}\right\}.
\]
If $p\in I^{<\beta+\varepsilon}$ and
$I(\Phi_1(p))\geq\beta-\varepsilon$, the entire trajectory between times
zero and one remains in the band
$|I-\beta|\leq\varepsilon<a/2$, where $\chi=1$. For
\[
 A_p:=\{t\in[0,1]\mid
 \Phi_t(p)\notin N_{\beta,\delta}\}
\]
we have, using $\delta\leq1$,
\[
 \frac{d}{dt}I(\Phi_t(p))\leq-\delta^2
 \quad(t\in A_p),
 \qquad
 I(\Phi_1(p))<\beta+\varepsilon-\delta^2|A_p|.
\]

Suppose that $p\notin N$ or $\Phi_1(p)\notin N$.
If $K_\beta=\varnothing$, choose $\delta$ with
$N_{\beta,\delta}=\varnothing$ and $A_p=[0,1]$. Otherwise,
set $h(t):=d(\Phi_t(p),K_\beta)$. By
\eqref{eq:velocidad-flujo-deformacion-finsler},
$|h(t)-h(s)|\leq2|t-s|$. Thus,
\[
 \begin{aligned}
 p\notin N&\ \Longrightarrow\
 h(t)\geq2\rho-2t\geq\rho
 &&(0\leq t\leq\rho/2),\\
 \Phi_1(p)\notin N&\ \Longrightarrow\
 h(t)\geq2\rho-2(1-t)\geq\rho
 &&(1-\rho/2\leq t\leq1).
 \end{aligned}
\]
Since $N_{\beta,\delta}\subseteq U_{\beta,\rho}$ and $\rho\leq1$,
either condition implies $|A_p|\geq\rho/2$.
Therefore,
\[
 I(\Phi_1(p))
 <\beta+\varepsilon-\frac{\rho\delta^2}{2}
 <\beta-\varepsilon.
\]
Applying this contradiction first when $p\notin N$ and then when
$\Phi_1(p)\notin N$ gives the two inclusions. The fixed-point properties
follow from uniqueness of the flow at zeros of $\mathbf{e}$.
\end{proof}

This theorem is the version of Theorem~II.3.11 in
\cite[p.~87]{Struwe2008}. The preceding proof displays its three global
ingredients separately: Palais--Smale sequential compactness, the
locally finite partition producing the pseudogradient, and Finsler completeness,
which prevents trajectories from escaping in finite time.

\section{Motivation and the mountain pass theorem}
\label{sec:preparacion-minimax-paso-montana}

The direct method seeks the bottom of the energy. In many superlinear
problems, however, the functional is not bounded below: along
certain rays the potential part eventually dominates and the
energy tends to $-\infty$. At the same time, the origin may be a
strict local minimum. To travel from this local valley to a region of
negative energy, every path must climb over a barrier. The mountain
pass selects, among all such paths, the smallest possible maximum height.
It does not minimize $I$ over $X$; it minimizes the maximum of $I$ along each path.

\begin{semblanzaHistorica}{Ambrosetti and Rabinowitz: a critical point by minimax}
The mountain pass theorem, introduced by Ambrosetti and Rabinowitz,
turned this geometric picture into an existence principle in infinite
dimensions~\cite{AmbrosettiRabinowitz1973}. Topology forces paths to
cross a barrier; the deformation lemma shows that the lowest crossing
height cannot be regular; and Palais--Smale prevents almost critical
sequences from escaping. The three ingredients are distinct, and each plays a
precise role.
\end{semblanzaHistorica}

\subsection{Paths, barriers, and the minimax level}

\begin{definition}[Path minimax level]
\label{def:nivel-minimax-caminos}
\index{minimax}
\index{mountain pass}
Let $u_0,u_1\in X$. The class of paths with fixed endpoints is
\[
 \Gamma(u_0,u_1):=
 \{\gamma\in C([0,1],X)\mid
 \gamma(0)=u_0,\ \gamma(1)=u_1\}.
\]
For $I\in C^1(X,\mathbb R)$, its minimax level is
\begin{equation}
 c(u_0,u_1):=
 \inf_{\gamma\in\Gamma(u_0,u_1)}
 \max_{t\in[0,1]}I(\gamma(t)).
\label{eq:nivel-minimax-caminos}
\end{equation}
\end{definition}

The class is nonempty because it contains the segment
$t\mapsto(1-t)u_0+tu_1$. For each path, the maximum of
$I\circ\gamma$ exists: it is the maximum of a continuous function on the
compact set $[0,1]$. Moreover, all these maxima are at least
$\displaystyle\max\{I(u_0),I(u_1)\}$, and the segment between the endpoints provides a
finite upper bound for their infimum. Thus the number in
\eqref{eq:nivel-minimax-caminos} is real.

\begin{lemma}[Forced crossing of the barrier]
\label{lem:cruce-barrera-paso-montana}
Let $u_0,u_1\in X$ and $r>0$, with $\|u_1-u_0\|>r$. Every
$\gamma\in\Gamma(u_0,u_1)$ intersects the sphere
\[
 S_r(u_0):=\{u\in X\mid\|u-u_0\|=r\}.
\]
Consequently, if
\[
 b:=\inf_{u\in S_r(u_0)}I(u),
\]
then $c(u_0,u_1)\geq b$.
\end{lemma}

\begin{proof}
The function $h(t):=\|\gamma(t)-u_0\|$ is continuous,
$h(0)=0$, and $h(1)=\|u_1-u_0\|>r$. The intermediate value theorem gives
$t_r\in(0,1)$ with $h(t_r)=r$. Hence,
\[
 \max_{t\in[0,1]}I(\gamma(t))
 \geq I(\gamma(t_r))
 \geq b.
\]
Taking the infimum over $\gamma$ proves the last assertion.
\end{proof}

\subsection{The abstract theorem and its proof}

\begin{theorem}[Mountain pass theorem]
\label{teo:paso-montana-banach-preparacion}
Let $X$ be a real Banach space, and let $I\in C^1(X,\mathbb R)$. Suppose
there exist $u_0,u_1\in X$ and $r>0$ such that
\[
 \|u_1-u_0\|>r
\]
and
\begin{equation}
 \max\{I(u_0),I(u_1)\}
 <b:=\inf_{\|u-u_0\|=r}I(u).
\label{eq:geometria-paso-montana}
\end{equation}
Let
\[
 c:=c(u_0,u_1)
 =\inf_{\gamma\in\Gamma(u_0,u_1)}
 \max_{t\in[0,1]}I(\gamma(t)).
\]
Then $c\geq b$ and, if $I$ satisfies $(PS)_c$, there exists $u_c\in X$ such that
\[
 I(u_c)=c,
 \qquad
 DI(u_c)=0.
\]
\end{theorem}

\begin{proof}
Lemma~\ref{lem:cruce-barrera-paso-montana} gives $c\geq b$. On the other hand,
the path class contains the segment between the endpoints, on which the
energy has a finite maximum. Therefore, $c\in\mathbb R$. Moreover,
\begin{equation}
 d:=c-\max\{I(u_0),I(u_1)\}>0.
\label{eq:separacion-extremos-nivel-minimax}
\end{equation}

Suppose, for a contradiction, that $K_c=\varnothing$. Since $I$ satisfies
$(PS)_c$, Proposition~\ref{prop:relacion-lemas-deformacion} provides
$a_0,\sigma>0$ such that
\[
 |I(u)-c|\leq2a_0
 \quad\Longrightarrow\quad
 \|DI(u)\|_{X'}\geq\sigma.
\]
Shrinking $a_0$ preserves this property. Take
\[
 0<a<\min\left\{a_0,\frac d2\right\}
 \qquad\text{and}\qquad
 0<\varepsilon<a.
\]
Theorem~\ref{teo:lema-deformacion-cuantitativo-banach}, applied with
$\beta=c$, produces a deformation $H$. By
\eqref{eq:separacion-extremos-nivel-minimax},
\[
 |I(u_i)-c|=c-I(u_i)\geq d>2a
 \qquad(i=0,1),
\]
so $H_1$ fixes both endpoints.

The definition of the infimum ensures that there exists
$\gamma_\varepsilon\in\Gamma(u_0,u_1)$ such that
\begin{equation}
 \max_{t\in[0,1]}I(\gamma_\varepsilon(t))<c+\varepsilon.
\label{eq:camino-casi-optimo-montana}
\end{equation}
Define
\[
 \widetilde\gamma_\varepsilon(t):=
 H_1(\gamma_\varepsilon(t)).
\]
Continuity of $H_1$ and the fact that the endpoints are fixed show that
$\widetilde\gamma_\varepsilon\in\Gamma(u_0,u_1)$. Moreover,
\eqref{eq:camino-casi-optimo-montana} and the conclusion of the quantitative lemma
give
\[
 I(\widetilde\gamma_\varepsilon(t))
 \leq c-\varepsilon
 \qquad(0\leq t\leq1).
\]
Consequently,
\[
 c
 \leq\max_{t\in[0,1]}I(\widetilde\gamma_\varepsilon(t))
 \leq c-\varepsilon,
\]
which is impossible. Therefore, $K_c\neq\varnothing$ and any of its
elements is the desired point $u_c$.
\end{proof}

The proof displays the role of each hypothesis. The topological barrier puts
$c$ above the endpoint energies; $(PS)_c$ turns the
absence of critical points into a uniform slope bound; and the deformation
lemma lowers an almost optimal path without moving its endpoints,
contradicting the definition of $c$. One may also directly apply
Struwe's lemma~\ref{teo:lema-deformacion-struwe-banach} with
$N=\varnothing$; Proposition~\ref{prop:relacion-lemas-deformacion} explains
why both arguments use the same mechanism in this case.

\begin{corollary}[Usual form with the origin and negative energy]
\label{cor:paso-montana-forma-usual}
Let $I\in C^1(X,\mathbb R)$ and suppose that:
\begin{enumerate}[label=(\roman*)]
\item $I(0)=0$;
\item there exist $\rho,\alpha>0$ such that
$I(u)\geq\alpha$ when $\|u\|=\rho$;
\item there exists $e\in X$ such that $\|e\|>\rho$ and $I(e)<0$.
\end{enumerate}
If
\[
 c_e:=\inf_{\gamma\in\Gamma(0,e)}
 \max_{t\in[0,1]}I(\gamma(t))
\]
and $I$ satisfies $(PS)_{c_e}$, then $c_e\geq\alpha$ and there exists
$u\in X\setminus\{0\}$ with
\[
 DI(u)=0,
 \qquad I(u)=c_e.
\]
\end{corollary}

\begin{proof}
Apply the theorem with $u_0=0$, $u_1=e$, $r=\rho$, and $b\geq\alpha$.
Since $c_e\geq\alpha>0=I(0)$, the resulting critical point cannot be the
origin.
\end{proof}

\subsection{Model: the Beltrami \texorpdfstring{$p$}{p}-Laplacian with a pure power}
\label{sec:modelo-potencia-paso-montana}

Let $(M,\mathbf g)$ be a nonempty compact Riemannian manifold of dimension
$n$, with smooth boundary, such that every connected component has nonempty boundary. Suppose
\[
 1<p<n,
 \qquad
 p<r<p^*:=\frac{np}{n-p},
\]
and equip
\[
 X:=W_0^{1,p}(M)
\]
with the equivalent norm
\[
 \|u\|:=\left(\int_M|\nabla u|_{\mathbf g}^p
 \,d\lambda_{\mathbf g}\right)^{1/p}.
\]
Equivalence with the Sobolev norm follows from the Poincaré inequality,
Theorem~\ref{teo:poincare-W0-haces}, and holds for every
$1<p<\infty$. The condition $p<n$ is imposed here to define the
finite critical exponent $p^*$. Consider
\begin{equation}
 I(u):=\frac1p\|u\|^p
 -\frac1r\|u\|_{L^r(M)}^r.
\label{eq:energia-potencia-p-paso-montana}
\end{equation}
Corollary~\ref{cor:potencia-norma-Lp-C1-haces}, applied to $TM$ with exponent $p$ and to the trivial real bundle with exponent $r$, together with the continuous linear operators $u\mapsto\nabla u$ and $X\hookrightarrow L^r(M)$, allows the chain rule to be applied. Thus, the functional belongs to $C^1(X,\mathbb R)$ and
\begin{equation}
 DI(u)[v]
 =\int_M|\nabla u|_{\mathbf g}^{p-2}
 \langle\nabla u,\nabla v\rangle_{\mathbf g}
 \,d\lambda_{\mathbf g}
 -\int_M|u|^{r-2}uv\,d\lambda_{\mathbf g}.
\label{eq:derivada-potencia-p-paso-montana}
\end{equation}
Its critical points are therefore exactly the weak solutions of
\begin{equation}
 \begin{cases}
 -\operatorname{div}_{\mathbf g}
 \bigl(|\nabla u|_{\mathbf g}^{p-2}\nabla u\bigr)
 =|u|^{r-2}u&\text{in }M,\\
 u=0&\text{on }\partial M.
 \end{cases}
\label{eq:problema-potencia-p-paso-montana}
\end{equation}

\begin{proposition}[Mountain pass geometry]
\label{prop:geometria-potencia-p-paso-montana}
The functional in \eqref{eq:energia-potencia-p-paso-montana} satisfies the
geometric hypotheses of
Corollary~\ref{cor:paso-montana-forma-usual}.
\end{proposition}

\begin{proof}
The continuous embedding $X\hookrightarrow L^r(M)$ provides $C_S>0$ such that
$\|u\|_{L^r(M)}\leq C_S\|u\|$. Therefore,
\begin{equation}
 I(u)\geq\frac1p\|u\|^p
 -\frac{C_S^r}{r}\|u\|^r.
\label{eq:cota-inferior-geometria-potencia}
\end{equation}
Since $r>p$, we can choose $\rho>0$ so that
\[
 \frac{C_S^r}{r}\rho^{r-p}\leq\frac1{2p}.
\]
If $\|u\|=\rho$, the preceding estimate gives
\[
 I(u)\geq\frac{\rho^p}{2p}=:\alpha>0.
\]

On the other hand, fix $w\in X\setminus\{0\}$. For $t>0$,
\[
 I(tw)=\frac{t^p}{p}\|w\|^p
 -\frac{t^r}{r}\|w\|_{L^r(M)}^r.
\]
The second coefficient is positive and $r>p$; dividing by $t^r$ shows that
$I(tw)/t^r\to-\|w\|_{L^r(M)}^r/r$. Consequently,
$I(tw)\to-\infty$. For sufficiently large $t$,
$e:=tw$ satisfies $\|e\|>\rho$ and $I(e)<0$.
\end{proof}

The strong convergence required by the Palais--Smale condition will
follow from the next monotonicity property.

\begin{lemma}[Property of type \texorpdfstring{$(S_+)$}{(S+)}]
\label{lem:S-mas-p-laplaciano-modelo}
Define $\mathcal A_p\colon X\to X'$ by
\[
 \mathcal A_p(u)[v]
 :=\int_M|\nabla u|_{\mathbf g}^{p-2}
 \langle\nabla u,\nabla v\rangle_{\mathbf g}
 \,d\lambda_{\mathbf g}.
\]
If $u_j\rightharpoonup u$ in $X$ and
\[
 \limsup_{j\to\infty}
 \mathcal A_p(u_j)[u_j-u]\leq0,
\]
then $u_j\to u$ in $X$.
\end{lemma}

\begin{proof}
Let $F(v):=\|v\|^p/p$. This functional is convex, is of class $C^1$, and
$DF(v)=\mathcal A_p(v)$. The supporting inequality for a differentiable convex
functional, applied at $u_j$, gives
\[
 F(u)\geq F(u_j)+DF(u_j)[u-u_j].
\]
Therefore,
\[
 \limsup_{j\to\infty}F(u_j)\leq F(u).
\]
Weak lower semicontinuity of the norm gives the opposite inequality for
the lower limit. It follows that
$\|u_j\|\to\|u\|$.

With the chosen norm, $X$ is identified isometrically through
$u\mapsto\nabla u$ with a closed subspace of
$L^p(M,TM)$. The latter is uniformly convex for
$1<p<\infty$, and uniform convexity passes to subspaces. Recall
why, in a uniformly convex space, weak convergence together with
convergence of norms implies strong convergence. If $u=0$, the assertion is
immediate. If $u\neq0$, normalize
$y_j:=u_j/\|u_j\|$ and $y:=u/\|u\|$. If a subsequence satisfied
$\|y_j-y\|\geq\eta>0$, uniform convexity would give
$\|(y_j+y)/2\|\leq1-\delta$ for some $\delta>0$. By Hahn--Banach there exists
$\ell\in X'$ with $\|\ell\|=1$ and $\ell(y)=1$; since $y_j\rightharpoonup y$,
\[
 \left\|\frac{y_j+y}{2}\right\|
 \geq\ell\left(\frac{y_j+y}{2}\right)
 \longrightarrow1,
\]
a contradiction. Thus, $y_j\to y$ and, using
$\|u_j\|\to\|u\|$, we obtain $u_j\to u$.
\end{proof}

\begin{proposition}[Palais--Smale for the subcritical power]
\label{prop:PS-potencia-p-paso-montana}
The functional in \eqref{eq:energia-potencia-p-paso-montana} satisfies
condition (PS).
\end{proposition}

\begin{proof}
Let $(u_j)\subset X$ be a Palais--Smale sequence. There exist $C>0$ and a
sequence $d_j\to 0^{+}$ such that
\[
 |I(u_j)|\leq C,
 \qquad
 \|DI(u_j)\|_{X'}\leq d_j.
\]
The exact identity
\begin{align}
 I(u_j)-\frac1rDI(u_j)[u_j]
 &=\left(\frac1p-\frac1r\right)\|u_j\|^p
\label{eq:identidad-acotacion-PS-potencia}
\end{align}
implies
\[
 \left(\frac1p-\frac1r\right)\|u_j\|^p
 \leq C+\frac{d_j}{r}\|u_j\|.
\]
The sequence is bounded. Indeed, if a subsequence had norm tending
to infinity, dividing the last inequality by $\|u_j\|^p$ would leave the
left-hand side positive while the right-hand side tended to zero, since $p>1$.

Since $X$ is reflexive, after passing to a subsequence there exists $u\in X$
such that $u_j\rightharpoonup u$. Corollary~\ref{cor:encajes-W0-haces-frontera} gives the compact embedding
because $r<p^*$; therefore,
\[
 u_j\longrightarrow u\quad\text{in }L^r(M).
\]
From \eqref{eq:derivada-potencia-p-paso-montana} we obtain
\begin{align*}
 \mathcal A_p(u_j)[u_j-u]
 &=DI(u_j)[u_j-u]
 +\int_M|u_j|^{r-2}u_j(u_j-u)\,d\lambda_{\mathbf g}.
\end{align*}
The first term tends to zero because $(u_j-u)$ is bounded in $X$ and
$\|DI(u_j)\|_{X'}\to0$. For the second, Hölder gives
\[
 \left|\int_M|u_j|^{r-2}u_j(u_j-u)\,d\lambda_{\mathbf g}\right|
 \leq\|u_j\|_{L^r(M)}^{r-1}
 \|u_j-u\|_{L^r(M)}\longrightarrow0.
\]
Lemma~\ref{lem:S-mas-p-laplaciano-modelo} implies $u_j\to u$ in $X$.
Every Palais--Smale sequence therefore has a convergent subsequence.
\end{proof}

\begin{corollary}[Nontrivial solution of the power problem]
\label{cor:solucion-potencia-p-paso-montana}
Problem \eqref{eq:problema-potencia-p-paso-montana} has a nontrivial
weak solution $u\in W_0^{1,p}(M)$.
\end{corollary}

\begin{proof}
Propositions~\ref{prop:geometria-potencia-p-paso-montana} and
\ref{prop:PS-potencia-p-paso-montana} allow
Corollary~\ref{cor:paso-montana-forma-usual} to be applied. The resulting critical point has
energy at least $\alpha>0$, so it is not the origin, and
\eqref{eq:derivada-potencia-p-paso-montana} is precisely the weak
formulation of \eqref{eq:problema-potencia-p-paso-montana}.
\end{proof}

This example clarifies when the mountain pass contributes something that
Nehari reduction does not automatically provide. The minimax argument uses only a
barrier, a low-energy direction, and compactness at the selected level;
it does not require every radial function $t\mapsto I(tu)$ to have a unique maximum or
$\{u\neq0\mid DI(u)[u]=0\}$ to be a submanifold. In the next chapter
we will compare the two methods and prove that, for the pure power, their
levels agree exactly.

The same abstract argument works on a complete Banach--Finsler
manifold whenever the minimax class is preserved by composition of its elements
with the deformation in
Theorem~\ref{teo:lema-deformacion-struwe-finsler}. This is the starting
point for more general geometric minimax principles.

\chapter{Variational constraints and the Nehari method}
\label{cap:restricciones-variacionales-nehari}

Many energy functionals cannot be minimized directly on the entire ambient space: the origin may be a trivial minimum, the functional may lack coercivity, or the relevant geometry may be concentrated on a constraint set. In such situations it is useful to study critical points of the functional on a suitable submanifold.

The structure of regular level sets and the description of their tangent spaces were established in Chapter~\ref{cap:variedades-banach}. Corollary~\ref{teo:multiplicadores-lagrange-banach}, obtained from Theorem~\ref{teo:funcion-implicita-banach}, provides the Lagrange multiplier principle in Banach spaces.\glsadd{multiplicador-lagrange} Here we apply that framework to variational constraints and, in particular, to the Nehari manifold.

In the Nehari method, the constraint is obtained by setting the derivative of the
functional in the radial direction $u$ equal to zero. This condition eliminates the direction in
which the energy may lose coercivity or lead to the trivial critical
point, while retaining the tangential variations relevant to locating
nonzero solutions.

\begin{semblanzaHistorica}{Zeev Nehari and natural constraints}
The method associated with Nehari uses the radial derivative of the energy to select a set where an unproductive direction disappears~\cite{Nehari1960}. Under suitable hypotheses, this constraint is natural: a critical point of the restricted functional is also critical on the whole space. The advantage is not merely technical. The radial identity reflects the balance between terms of different homogeneities and turns the search for nontrivial solutions into a geometric problem on a submanifold of the function space.
\end{semblanzaHistorica}

\section{The Nehari method}
\label{sec:metodo-nehari}
\index{method!Nehari}

The Nehari manifold turns a necessary variational condition into a
useful geometric constraint. Let $I\colon X\longrightarrow \mathbb{R}$ be an energy functional
on a Banach space suitable for seeking weak solutions.
Every critical point $u$ satisfies $DI(u)=0$ in $X'$ and, in particular,
$DI(u)(u)=0$.

This motivates considering the set
\[
\mathcal{N}:=\{u\in X\setminus \{0\}\mid DI(u)(u)=0\},
\]
called the \emph{Nehari manifold associated with $I$}.\index{manifold!Nehari}\glsadd{variedad-nehari} The idea of the method is
to study first the functional $I$ restricted to $\mathcal{N}$ and then
to show that, under suitable hypotheses, critical points of
$I\restriction_{\mathcal{N}}$ are in fact critical points of $I$ on all of $X$.

The results of Section~\ref{sec:conjuntos-nivel-regulares-banach} provide a geometric description of this constraint. If $V=\{u\in\Omega\mid J(u)=c\}$ is a regular level set in an open set $\Omega\subseteq X$, then $T_uV=\operatorname{ker}(DJ(u))$ and
\[
d(\varphi\restriction_V)_u
=
D\varphi(u)\restriction_{T_uV}.
\]
For this reason, define
\[
\|D\varphi(u)\|_{*,V}
:=
\|D\varphi(u)\restriction_{T_uV}\|_{\mathcal L(T_uV,\mathbb{R})}.
\]
If $T_uV\neq\{0\}$, this norm can also be written as
\[
\|D\varphi(u)\|_{*,V}
=
\sup_{\substack{v\in T_uV\\ \|v\|=1}}|D\varphi(u)(v)|.
\]
The norm on $T_uV$ is induced by its identification with the subspace $\operatorname{ker}(DJ(u))\subseteq X$ and measures the slope of $\varphi$ only in admissible directions.

\begin{lemma}\label{lem:funcionales-nucleo-incluido}
Let $X$ be a real vector space, and let $f,g\colon X\longrightarrow \mathbb{R}$ be
linear functionals, with $g\neq 0$. If
$\operatorname{ker}(g)\subset\operatorname{ker}(f)$, then there exists a unique
$\lambda\in\mathbb{R}$ such that $f=\lambda g$.
\end{lemma}

\begin{proof}
Since $g\neq 0$, there exists $x_{0}\in X$ such that $g(x_{0})\neq 0$. Define
$\lambda:=\displaystyle\frac{f(x_{0})}{g(x_{0})}$. We show that $f(x)=\lambda g(x)$ for every
$x\in X$. Let $x\in X$ and consider
$y:=x-\displaystyle\frac{g(x)}{g(x_{0})}x_{0}$. Then
$g(y)=g(x)-\displaystyle\frac{g(x)}{g(x_{0})}g(x_{0})=0$, so
$y\in\operatorname{ker}(g)$. By hypothesis, $y\in\operatorname{ker}(f)$, and
then
$0=f(y)=f(x)-\displaystyle\frac{g(x)}{g(x_{0})}f(x_{0})$. Therefore,
$f(x)=\displaystyle\frac{f(x_{0})}{g(x_{0})}g(x)=\lambda g(x)$. Since this holds for every
$x\in X$, it follows that $f=\lambda g$. Uniqueness is immediate: if
$\lambda_{1}g=\lambda_{2}g$, then $(\lambda_{1}-\lambda_{2})g=0$, and since
$g\neq 0$, it follows that $\lambda_{1}=\lambda_{2}$.
\end{proof}

\begin{lemma}[Duality lemma]\label{lem:dualidad-norma-restringida}\index{duality lemma}
Let $X$ be a real Banach space, and let $f,g\in X'$, with $g\neq0$. Then
\[
\|f\restriction_{\operatorname{ker}(g)}\|_{(\operatorname{ker}(g))'}
=
\min_{\lambda\in\mathbb{R}}\|f-\lambda g\|_{X'}.
\]
If $\operatorname{ker}(g)\neq\{0\}$, the left-hand side agrees with $\displaystyle\sup_{\substack{y\in\operatorname{ker}(g)\\ \|y\|=1}}|f(y)|$.
\end{lemma}

\begin{proof}
Let $\lambda\in\mathbb{R}$. Since $(f-\lambda g)\restriction_{\operatorname{ker}(g)}=f\restriction_{\operatorname{ker}(g)}$, we have
\[
\|f\restriction_{\operatorname{ker}(g)}\|_{(\operatorname{ker}(g))'}
\leq
\|f-\lambda g\|_{X'}.
\]
Taking the infimum over $\lambda$, we obtain
\[
\|f\restriction_{\operatorname{ker}(g)}\|_{(\operatorname{ker}(g))'}
\leq
\inf_{\lambda\in\mathbb{R}}\|f-\lambda g\|_{X'}.
\]

By Corollary~\ref{cor:extension-hahn-banach-preserva-norma}, there exists $\widetilde f\in X'$ extending $f\restriction_{\operatorname{ker}(g)}$ and satisfying
\[
\|\widetilde f\|_{X'}
=
\|f\restriction_{\operatorname{ker}(g)}\|_{(\operatorname{ker}(g))'}.
\]
Since $f-\widetilde f$ vanishes on $\operatorname{ker}(g)$, Lemma~\ref{lem:funcionales-nucleo-incluido} provides $\lambda_{0}\in\mathbb{R}$ such that $f-\widetilde f=\lambda_{0}g$. Therefore,
\[
\inf_{\lambda\in\mathbb{R}}\|f-\lambda g\|_{X'}
\leq
\|f-\lambda_{0}g\|_{X'}
=
\|\widetilde f\|_{X'}
=
\|f\restriction_{\operatorname{ker}(g)}\|_{(\operatorname{ker}(g))'}.
\]
The two inequalities prove equality and also show that the infimum is attained at $\lambda_{0}$. The last assertion is the usual characterization of the norm of a continuous linear functional on a nontrivial normed space.
\end{proof}

\begin{proposition}[Characterization of constrained critical points]\label{prop:puntos-criticos-restricciones}\index{characterization of constrained critical points@characterization of constrained critical points}
Let $X$ be a real Banach space, let $\Omega\subset X$ be an open set,
and let $J\in C^{1}(\Omega,\mathbb{R})$ and $c\in\mathbb{R}$. Suppose that
$V:=\{u\in\Omega\mid J(u)=c\}$ is a regular level set, that is,
$DJ(u)\neq 0$ for every $u\in V$. If
$\varphi\in C^{1}(\Omega,\mathbb{R})$ and $u\in V$, then

\[
\|D\varphi(u)\|_{*,V}
=
\min_{\lambda\in\mathbb{R}}
\|D\varphi(u)-\lambda DJ(u)\|_{X'}.
\]

Moreover, $u$ is a critical point of the restriction $\varphi\restriction_{V}$ if and only if
there exists $\lambda\in\mathbb{R}$ such that

\[
D\varphi(u)=\lambda DJ(u).
\]
\end{proposition}

\begin{proof}
Since $V$ is a regular level set of $J$ in the open set $\Omega$, Theorem~\ref{teo:conjuntos-nivel-regulares-banach} implies that $T_uV=\operatorname{ker}(DJ(u))$. Moreover, by Lemma~\ref{lem:diferencial-restriccion-subvariedad},
\[
d(\varphi\restriction_V)_u
=
D\varphi(u)\restriction_{T_uV}.
\]
Therefore,
\[
\|D\varphi(u)\|_{*,V}
=
\|D\varphi(u)\restriction_{\operatorname{ker}(DJ(u))}\|_{(\operatorname{ker}(DJ(u)))'}.
\]
Applying Lemma~\ref{lem:dualidad-norma-restringida} with $f=D\varphi(u)$ and $g=DJ(u)$ gives

\[
\|D\varphi(u)\|_{*,V}
=
\min_{\lambda\in\mathbb{R}}\|D\varphi(u)-\lambda DJ(u)\|_{X'}.
\]

We now prove the characterization of critical points. By definition, $u$ is a
critical point of $\varphi\restriction_{V}$ if and only if $d(\varphi\restriction_{V})_{u}=0$. By
Lemma~\ref{lem:diferencial-restriccion-subvariedad}, this is equivalent to
$D\varphi(u)(v)=0$ for every $v\in T_{u}V$. Since
$T_{u}V=\operatorname{ker}(DJ(u))$, this is equivalent to
$\operatorname{ker}(DJ(u))\subset \operatorname{ker}(D\varphi(u))$. Since
$DJ(u)\neq 0$, Lemma~\ref{lem:funcionales-nucleo-incluido} gives
$\lambda\in\mathbb{R}$ such that $D\varphi(u)=\lambda DJ(u)$.

Conversely, if there exists $\lambda\in\mathbb{R}$ such that
$D\varphi(u)=\lambda DJ(u)$, then for every
$v\in T_{u}V=\operatorname{ker}(DJ(u))$ we have
$D\varphi(u)(v)=\lambda DJ(u)(v)=0$. Therefore,
$d(\varphi\restriction_{V})_{u}=0$ and $u$ is a critical point of $\varphi\restriction_{V}$. This proves
the equivalence.
\end{proof}

In particular, for the Nehari manifold take
$\Omega:=X\setminus\{0\}$, which is open in $X$, and define
$J\colon \Omega\longrightarrow\mathbb R$ by $J(u):=DI(u)(u)$. Then

\[
\mathcal N
=
\{u\in\Omega\mid J(u)=0\}
=
\{u\in X\setminus\{0\}\mid DI(u)(u)=0\}.
\]

Therefore, if $J\in C^1(\Omega,\mathbb R)$ and $DJ(u)\neq0$ for every
$u\in\mathcal N$, Theorem~\ref{teo:conjuntos-nivel-regulares-banach} implies that $\mathcal N$ is an
embedded Banach submanifold of class $C^1$ and codimension $1$. Furthermore, for each
$u\in\mathcal N$ we have

\[
T_u\mathcal N=\operatorname{ker}(DJ(u)).
\]

If $I$ is also of class $C^2$ on $\Omega$, then $J$ is of class $C^1$ and,
for every $v\in X$,

\[
DJ(u)(v)=D^2I(u)(v,u)+DI(u)(v).
\]

In the next section we will apply this framework to the superlinear subcritical problem on a compact Riemannian manifold.

\section{The Nehari manifold in a model case}
\label{sec:nehari-caso-modelo}

Let $(M,\mathbf{g})$ be a nonempty compact Riemannian manifold with smooth boundary. Consider $X=W_0^{1,2}(M)$ and
\[
I(u):=\frac{1}{2}\int_M|\nabla u|_{\mathbf{g}}^2\,d\lambda_{\mathbf{g}}-\frac{1}{r}\int_M|u|^r\,d\lambda_{\mathbf{g}},
\qquad 2<r<2^*.
\]
Then
\[
J(u):=DI(u)(u)=\int_M|\nabla u|_{\mathbf{g}}^2\,d\lambda_{\mathbf{g}}-\int_M|u|^r\,d\lambda_{\mathbf{g}},
\]
and the Nehari manifold is $\mathcal N=J^{-1}(0)\setminus\{0\}$.

\begin{proposition}[The Nehari manifold as a hypersurface]
\label{prop:nehari-variedad-caso-modelo}
Let $(M,\mathbf{g})$ be a compact Riemannian manifold with nonempty boundary, and let $X$, $I$, $J$, and $\mathcal N$ be the objects defined at the beginning of this section. The set $\mathcal N$ is an embedded Banach submanifold of class $C^1$ and codimension $1$ in $X\setminus\{0\}$. Moreover, for each $u\in\mathcal N$,
\[
T_u\mathcal N=\ker DJ(u).
\]
\end{proposition}

\begin{proof}
By Proposition~\ref{prop:funcional-energia-escalar-C2}, $I\in C^2(X,\mathbb R)$ and $J\in C^1(X\setminus\{0\},\mathbb R)$. For $u,v\in X$,
\[
DJ(u)(v)=2\int_M\langle\nabla u,\nabla v\rangle_{\mathbf{g}}\,d\lambda_{\mathbf{g}}-r\int_M|u|^{r-2}uv\,d\lambda_{\mathbf{g}}.
\]
If $u\in\mathcal N$, then
\[
DJ(u)(u)=(2-r)\int_M|u|^r\,d\lambda_{\mathbf{g}}<0.
\]
Therefore, $DJ(u)\neq0$ for every $u\in\mathcal N$. Theorem~\ref{teo:conjuntos-nivel-regulares-banach} implies that $\mathcal N$ is an embedded submanifold of codimension $1$ and that $T_u\mathcal N=\ker DJ(u)$.
\end{proof}

This example explains why the term Nehari manifold is used when the functional is of class $C^2$. In quasilinear problems, however, the functional may be only of class $C^1$; in that situation it is not legitimate to assume that $u\mapsto DI(u)(u)$ is differentiable. The radial method of Szulkin and Weth works directly with the topological structure of the constraint.

\section{The radial reduction of Szulkin and Weth}
\label{sec:reduccion-radial-szulkin-weth}

The Nehari method can be developed without assuming that the constraint set is a differentiable submanifold. The idea is to select on each positive ray the unique point where the functional attains its maximum and to transfer the variational problem to the unit sphere. We present the abstract framework of Szulkin and Weth~\cite[\S3.1, Propositions~8--9 and Corollary~10]{SzulkinWethNehari} in the notation adopted in this book.

Let $X$ be a uniformly convex real Banach space, let $I\in C^1(X,\mathbb R)$, and suppose that $I(0)=0$. Denote the unit sphere by
\[
S_X:=\{u\in X\mid\|u\|_X=1\}
\]
A function $\varphi\colon [0,+\infty)\longrightarrow[0,+\infty)$ is called a \textbf{normalization function} if it is continuous, $\varphi(0)=0$, it is strictly increasing, and $\varphi(t)\to+\infty$ as $t\to+\infty$. We consider the following hypotheses:
\begin{enumerate}[label=\textup{(A\arabic*)},ref=A\arabic*]
\item\label{hip:szulkin-A1} There exists a normalization function $\varphi$ such that the functional
\[
\psi(u):=\int_0^{\|u\|_X}\varphi(t)\,dt,
\qquad u\in X,
\]
belongs to $C^1(X\setminus\{0\},\mathbb R)$, the duality map $\mathcal J:=D\psi\colon X\setminus\{0\}\longrightarrow X'$ maps bounded sets to bounded sets, and
\[
\mathcal J(w)(w)=1
\]
for every $w\in S_X$.
\item\label{hip:szulkin-A2} For each $w\in X\setminus\{0\}$ there exists $s_w>0$ such that, if $\alpha_w\colon (0,+\infty)\longrightarrow\mathbb R$ is given by $\alpha_w(s):=I(sw)$, then
\[
\alpha_w'(s)>0\quad\text{if }0<s<s_w,
\qquad
\alpha_w'(s)<0\quad\text{if }s>s_w.
\]
\item\label{hip:szulkin-A3} There exists $\delta>0$ such that $s_w\geq\delta$ for every $w\in S_X$ and, for each compact set $K\subseteq S_X$, there exists $C_K>0$ such that $s_w\leq C_K$ for every $w\in K$.
\end{enumerate}

Hypothesis \textup{(\ref{hip:szulkin-A2})} implies that $s_w$ is the unique global maximum of $\alpha_w$ and that $s_ww$ is the unique point of the positive ray generated by $w$ belonging to the Nehari set
\[
\mathcal N:=\{u\in X\setminus\{0\}\mid DI(u)(u)=0\}.
\]

\begin{theorem}[Radial reduction principle]
\label{teo:reduccion-radial-szulkin-weth}
Suppose hypotheses \textup{(\ref{hip:szulkin-A1})}--\textup{(\ref{hip:szulkin-A3})} hold. Define
\[
\widehat m\colon X\setminus\{0\}\longrightarrow\mathcal N,
\qquad
\widehat m(w):=s_ww,
\]
and let $m:=\widehat m\restriction_{S_X}$. Then the following assertions hold:
\begin{enumerate}
\item the map $\widehat m$ is continuous;
\item $m\colon S_X\longrightarrow\mathcal N$ is a homeomorphism and
\[
m^{-1}(u)=\frac{u}{\|u\|_X}
\]
for every $u\in\mathcal N$;
\item the functional $\widehat\Psi:=I\circ\widehat m$ belongs to $C^1(X\setminus\{0\},\mathbb R)$ and
\begin{equation}
\label{eq:derivada-reduccion-radial-abierta}
D\widehat\Psi(w)(z)
=
\frac{\|\widehat m(w)\|_X}{\|w\|_X}
DI(\widehat m(w))(z)
\end{equation}
for all $w\in X\setminus\{0\}$ and $z\in X$;
\item the sphere $S_X$ is an embedded Banach submanifold of class $C^1$, the functional $\Psi:=I\circ m$ belongs to $C^1(S_X,\mathbb R)$, and
\begin{equation}
\label{eq:derivada-reduccion-radial}
D\Psi(w)(z)=\|m(w)\|_XDI(m(w))(z)
\end{equation}
for all $w\in S_X$ and $z\in T_wS_X$;
\item a point $w\in S_X$ is critical for $\Psi$ if and only if $m(w)$ is a nontrivial critical point of $I$;
\item the minimax characterization
\[
\inf_{u\in\mathcal N}I(u)
=
\inf_{w\in S_X}\Psi(w)
=
\inf_{w\in X\setminus\{0\}}\max_{s>0}I(sw).
\]
holds.
\end{enumerate}
\end{theorem}

\begin{proof}
Begin with a scaling property. If $w\in X\setminus\{0\}$ and $t>0$, then $\alpha_{tw}(s)=I(stw)=\alpha_w(st)$. Uniqueness of the maximum given by \textup{(\ref{hip:szulkin-A2})} implies
\begin{equation}
\label{eq:escala-parametro-nehari}
s_{tw}=\frac{s_w}{t}
\end{equation}
and hence $\widehat m(tw)=\widehat m(w)$.

The first part of \textup{(\ref{hip:szulkin-A3})} shows that $\mathcal N$ is bounded away from the origin. If $u\in\mathcal N$ and $w:=\frac{u}{\|u\|_X}\in S_X$, then $u=s_ww$ and $\|u\|_X=s_w\geq\delta$. Moreover, $\mathcal N$ is closed in $X$. If $u_j\in\mathcal N$ and $u_j\to u$ in $X$, separation from the origin implies $u\neq0$, whereas continuity of $DI\colon X\longrightarrow X'$ gives
\[
DI(u)(u)=\lim_{j\to\infty}DI(u_j)(u_j)=0.
\]
Consequently, $u\in\mathcal N$.

We prove continuity of $\widehat m$. Let $w_j\to w\neq0$. By \eqref{eq:escala-parametro-nehari}, we may replace each $w_j$ and $w$ by their normalizations and assume that $w_j,w\in S_X$. The set $K:=\{w\}\cup\{w_j\mid j\in\mathbb N\}$ is compact. By \textup{(\ref{hip:szulkin-A3})}, the sequence $(s_{w_j})$ is bounded and bounded away from zero. Every subsequence therefore has a further subsequence, denoted in the same way, such that $s_{w_j}\to s>0$. Hence,
\[
\widehat m(w_j)=s_{w_j}w_j\longrightarrow sw.
\]
Since $\mathcal N$ is closed, $sw\in\mathcal N$. Uniqueness of the intersection of the positive ray generated by $w$ with $\mathcal N$ implies $s=s_w$, and therefore $\widehat m(w_j)\to\widehat m(w)$. The subsequence criterion proves continuity of $\widehat m$.

Continuity of $m$ is immediate. If $u\in\mathcal N$ and $w:=\frac{u}{\|u\|_X}$, then $w\in S_X$ and, by radial uniqueness, $s_w=\|u\|_X$; thus, $m(w)=u$. This proves that $m$ is bijective and that its inverse is given by $m^{-1}(u)=\frac{u}{\|u\|_X}$. Since normalization is continuous on $X\setminus\{0\}$, $m^{-1}$ is continuous and $m$ is a homeomorphism.

We now study regularity of $\widehat\Psi$. Fix $w\in X\setminus\{0\}$ and $z\in X$. For sufficiently small $t\neq0$, write $s_t:=s_{w+tz}$. Since $s_ww$ maximizes $s\mapsto I(sw)$ and $s_t(w+tz)$ maximizes $s\mapsto I(s(w+tz))$, we have
\[
I(s_w(w+tz))-I(s_ww)
\leq
\widehat\Psi(w+tz)-\widehat\Psi(w)
\leq
I(s_t(w+tz))-I(s_tw).
\]
Applying the mean value theorem to the real functions obtained by restricting $I$ to the corresponding segments gives $\eta_t,\tau_t\in(0,1)$ such that
\[
I(s_w(w+tz))-I(s_ww)
=
t s_wDI\bigl(s_w(w+\eta_ttz)\bigr)(z)
\]
and
\[
I(s_t(w+tz))-I(s_tw)
=
t s_tDI\bigl(s_t(w+\tau_ttz)\bigr)(z).
\]
Continuity of $w\mapsto s_w$ and $DI$ implies that both quotients, after division by $t$, converge to $s_wDI(s_ww)(z)$ as $t\to0$. For $t>0$ the result follows directly by squeezing, and for $t<0$ the inequalities reverse on division by $t$ and produce the same limit. Therefore, the directional derivative exists and
\[
D_z\widehat\Psi(w)=s_wDI(s_ww)(z).
\]
Since $s_w=\frac{\|\widehat m(w)\|_X}{\|w\|_X}$, we obtain \eqref{eq:derivada-reduccion-radial-abierta}. The right-hand side is linear and continuous in $z$, and depends continuously on $w$ in the norm of $X'$. By Corollary~\ref{cor:caracterizacion-C1-gateaux}, $\widehat\Psi\in C^1(X\setminus\{0\},\mathbb R)$.

Let $c_\varphi:=\displaystyle\int_0^1\varphi(t)\,dt$. Since $\varphi$ is strictly increasing, $S_X=\psi^{-1}(c_\varphi)$. Moreover, \textup{(\ref{hip:szulkin-A1})} gives $D\psi(w)(w)=\mathcal J(w)(w)=1$ for every $w\in S_X$. The value $c_\varphi$ is regular, and Theorem~\ref{teo:conjuntos-nivel-regulares-banach} implies that $S_X$ is an embedded submanifold of class $C^1$ and
\[
T_wS_X=\ker\mathcal J(w).
\]
Identity \eqref{eq:derivada-reduccion-radial} follows by restricting \eqref{eq:derivada-reduccion-radial-abierta} to $S_X$, since $\|w\|_X=1$.

If $DI(m(w))=0$, then \eqref{eq:derivada-reduccion-radial} shows that $D\Psi(w)=0$. Conversely, suppose that $D\Psi(w)=0$. Given $v\in X$, define
\[
\alpha:=\mathcal J(w)(v),
\qquad
z:=v-\alpha w.
\]
Since $\mathcal J(w)(w)=1$, we have $\mathcal J(w)(z)=0$, so $z\in T_wS_X$. The Nehari identity gives $DI(m(w))(m(w))=0$ and, since $m(w)=\|m(w)\|_Xw$, also $DI(m(w))(w)=0$. On the other hand, \eqref{eq:derivada-reduccion-radial} and $D\Psi(w)=0$ imply $DI(m(w))(z)=0$. Consequently,
\[
DI(m(w))(v)=DI(m(w))(z)+\alpha DI(m(w))(w)=0.
\]
Since $v$ was arbitrary, $DI(m(w))=0$.

Finally, bijectivity of $m$ gives
\[
\inf_{u\in\mathcal N}I(u)=\inf_{w\in S_X}\Psi(w).
\]
By \textup{(\ref{hip:szulkin-A2})}, $\displaystyle I(s_ww)=\displaystyle\max_{s>0}I(sw)$, and this quantity depends only on the positive ray generated by $w$. Normalizing each nonzero vector then yields
\[
\inf_{w\in S_X}\Psi(w)
=
\inf_{w\in X\setminus\{0\}}\max_{s>0}I(sw).
\]
\end{proof}

Radial reduction allows minimization on $\mathcal N$, which need not be a differentiable submanifold, to be replaced by a problem of class $C^1$ on the unit sphere.

\section{Mountain pass and Nehari in the power model}
\label{sec:comparacion-montana-nehari-potencia}

Return to the power problem of
Chapter~\ref{cap:palais-smale-deformacion}. Let $(M,\mathbf g)$ be a nonempty compact
Riemannian manifold of dimension $n$, with smooth boundary, such that every connected component has nonempty boundary,
and suppose
\[
 1<p<n,
 \qquad
 p<r<p^*:=\frac{np}{n-p}.
\]
On
\[
 X:=W_0^{1,p}(M),
 \qquad
 \|u\|:=\|\nabla u\|_{L^p(M)},
\]
consider
\begin{equation}
 I(u):=\frac1p\|u\|^p-\frac1r\|u\|_{L^r(M)}^r.
\label{eq:energia-potencia-comparacion-nehari}
\end{equation}
Its derivative was calculated in
\eqref{eq:derivada-potencia-p-paso-montana}. Introduce
\[
 G(u):=DI(u)[u]
 =\|u\|^p-\|u\|_{L^r(M)}^r
\]
and
\[
 \mathcal N:=\{u\in X\setminus\{0\}\mid G(u)=0\}.
\]
The scaling argument allows each positive ray to be projected onto
$\mathcal N$ and relates this projection to energy minimization.

\begin{theorem}[Radial projection, minimizer, and ground state]
\label{teo:nehari-potencia-minimizador-completo}
For the functional in \eqref{eq:energia-potencia-comparacion-nehari}, the following
assertions hold.
\begin{enumerate}[label=(\alph*)]
\item For each $u\in X\setminus\{0\}$ there exists a unique $t_u>0$ such that
$t_uu\in\mathcal N$, given by
\begin{equation}
 t_u=\left(
 \frac{\|u\|^p}{\|u\|_{L^r(M)}^r}
 \right)^{\!1/(r-p)}.
\label{eq:escala-nehari-potencia-p}
\end{equation}
Moreover,
\[
 I(t_uu)=\max_{t>0}I(tu).
\]
\item The Nehari level
\[
 m:=\inf_{u\in\mathcal N}I(u)
\]
is positive and is attained at some $u_0\in\mathcal N$.
\item The minimizer $u_0$ is a nontrivial critical point of $I$ and is a
ground state:
\[
 I(u_0)=m
 =\inf\{I(v)\mid v\neq0,\ DI(v)=0\}.
\]
\item The radial characterization
\begin{equation}
 m=\inf_{u\in X\setminus\{0\}}\max_{t>0}I(tu).
\label{eq:nivel-nehari-minimax-radial-potencia}
\end{equation}
holds.
\end{enumerate}
\end{theorem}

\begin{proof}
Fix $u\neq0$ and define $\varphi_u(t):=I(tu)$ for $t>0$. Then
\begin{align*}
 \varphi_u(t)
 &=\frac{t^p}{p}\|u\|^p
 -\frac{t^r}{r}\|u\|_{L^r(M)}^r,\\
 \varphi_u'(t)
 &=t^{p-1}\left(
 \|u\|^p-t^{r-p}\|u\|_{L^r(M)}^r
 \right).
\end{align*}
Both coefficients are positive. Since $t\mapsto t^{r-p}$ is
strictly increasing, $\varphi_u'$ vanishes exactly once, precisely at
the value in \eqref{eq:escala-nehari-potencia-p}; it is positive before and negative
afterward. This point is therefore the unique global maximum of
$\varphi_u$. Moreover,
\[
 G(tu)=t\varphi_u'(t),
\]
which proves (a).

If $u\in\mathcal N$, the Nehari identity gives
\begin{equation}
 \|u\|^p=\|u\|_{L^r(M)}^r,
 \qquad
 I(u)=\left(\frac1p-\frac1r\right)\|u\|^p.
\label{eq:energia-sobre-nehari-potencia-p}
\end{equation}
Let $C_S>0$ be a constant for the embedding
$X\hookrightarrow L^r(M)$. Then
\[
 \|u\|^p=\|u\|_{L^r(M)}^r
 \leq C_S^r\|u\|^r.
\]
Since $u\neq0$, it follows that
\[
 \|u\|\geq C_S^{-r/(r-p)}=:\delta>0.
\]
Together with \eqref{eq:energia-sobre-nehari-potencia-p}, this shows that
\begin{equation}
 m\geq\left(\frac1p-\frac1r\right)\delta^p>0.
\label{eq:positividad-nivel-nehari-potencia-p}
\end{equation}

By (a), $\mathcal N$ is nonempty. For $v\in X\setminus\{0\}$,
we have $m\leq I(t_vv)<\infty$; combined with the preceding bound,
this gives $0<m<\infty$. For each $j\geq1$, choose $u_j\in\mathcal N$ with
\[
 m\leq I(u_j)<m+1/j.
\]
Identity \eqref{eq:energia-sobre-nehari-potencia-p} makes $(u_j)$ bounded
in $X$. Since $X$ is reflexive, a subsequence converges weakly to some
$u\in X$. The subcritical compact embedding gives a further subsequence
converging strongly to some $v\in L^r(M)$. The two limits coincide.
Indeed, for each $\psi\in L^{r'}(M)$, $r'=r/(r-1)$, the functional
$w\mapsto\int_Mw\psi\,d\lambda_{\mathbf g}$ is continuous on $X$ by
Hölder and the continuous embedding into $L^r(M)$. The two convergences give
\[
 \int_Mv\psi\,d\lambda_{\mathbf g}
 =\lim_{j\to\infty}\int_Mu_j\psi\,d\lambda_{\mathbf g}
 =\int_Mu\psi\,d\lambda_{\mathbf g}.
\]
Duality for $L^r(M)$ implies $v=u$ almost everywhere. Retaining the
notation for the final subsequence, we have
\[
 u_j\rightharpoonup u\quad\text{in }X,
 \qquad
 u_j\to u\quad\text{in }L^r(M).
\]
If $\kappa:=1/p-1/r>0$, then
\begin{equation}
 m=\lim_{j\to\infty}I(u_j)
 =\kappa\lim_{j\to\infty}\|u_j\|^p
 =\kappa\lim_{j\to\infty}\|u_j\|_{L^r(M)}^r
 =\kappa\|u\|_{L^r(M)}^r.
\label{eq:limite-minimizante-nehari-potencia}
\end{equation}
Positivity of $m$ implies $u\neq0$. Weak lower semicontinuity gives
\begin{equation}
 \|u\|^p
 \leq\liminf_{j\to\infty}\|u_j\|^p
 =\|u\|_{L^r(M)}^r.
\label{eq:desigualdad-limite-nehari-potencia}
\end{equation}

If the inequality were strict, formula
\eqref{eq:escala-nehari-potencia-p} would give $0<t_u<1$ and
$t_uu\in\mathcal N$. Using
\eqref{eq:energia-sobre-nehari-potencia-p} and
\eqref{eq:limite-minimizante-nehari-potencia} again, we would obtain
\[
 m\leq I(t_uu)
 =\kappa t_u^p\|u\|^p
 <\kappa\|u\|_{L^r(M)}^r
 =m,
\]
a contradiction. Thus equality holds in
\eqref{eq:desigualdad-limite-nehari-potencia}; consequently,
$u\in\mathcal N$ and $I(u)=m$. We have proved (b), with $u_0:=u$.

We now prove that the constraint is natural. The functional $G$ belongs to
$C^1(X,\mathbb R)$ and, for $v\in\mathcal N$,
\[
 DG(v)[v]
 =p\|v\|^p-r\|v\|_{L^r(M)}^r
 =(p-r)\|v\|^p<0.
\]
In particular, $DG(v)\neq0$ on $\mathcal N$, so $\mathcal N$ is
a hypersurface of class $C^1$. Since $u_0$ minimizes $I$ on this hypersurface, the
Lagrange multiplier principle provides $\lambda\in\mathbb R$
such that
\[
 DI(u_0)=\lambda DG(u_0).
\]
Evaluating at $u_0$ gives
\[
 0=DI(u_0)[u_0]
 =\lambda(p-r)\|u_0\|^p,
\]
and hence $\lambda=0$ and $DI(u_0)=0$. Every nontrivial critical point
$v$ of $I$ satisfies $G(v)=DI(v)[v]=0$, and therefore belongs to $\mathcal N$ with
$I(v)\geq m=I(u_0)$. This proves (c).

Finally, by (a), each positive ray meets $\mathcal N$ exactly
once and attains its maximum there. Taking all nonzero directions
gives
\[
 \inf_{u\neq0}\max_{t>0}I(tu)
 =\inf_{v\in\mathcal N}I(v)=m,
\]
which is (d).
\end{proof}

\subsection{Conceptual comparison}

The two methods organize the same energy in different ways. The following table
compares their hypotheses, constructions, and conclusions.

\begin{table}[htbp]
\centering
\caption{The mountain pass and Nehari methods.}
\label{tab:comparacion-montana-nehari}
\begingroup
\small
\renewcommand{\arraystretch}{1.18}
\begin{tabularx}{\textwidth}{>{\raggedright\arraybackslash}p{0.20\textwidth}
 >{\raggedright\arraybackslash}X
 >{\raggedright\arraybackslash}X}
\toprule
& \textbf{Mountain pass} & \textbf{Nehari method}\\
\midrule
Main idea
& Minimax over paths
& Minimization on a natural constraint\\
Main object
& The class $\Gamma(0,e)$
& The manifold $\mathcal N$\\
Variational level
& $\displaystyle c_e=\inf_{\gamma\in\Gamma(0,e)}\displaystyle\max_t I(\gamma(t))$
& $\displaystyle m=\inf_{u\in\mathcal N}I(u)$\\
Geometry
& Mountain pass barrier
& Radial geometry\\
Construction
& Paths from the origin to a point of negative energy
& Points satisfying $DI(u)[u]=0$\\
Compactness
& Usually $(PS)_{c_e}$
& Compactness of minimizing sequences\\
Typical result
& A nontrivial critical point
& Often a ground state\\
Main strength
& Flexible minimax principle; does not require radial uniqueness
& Precise description of least-energy solutions\\
Main difficulty
& Verifying Palais--Smale at the minimax level
& Proving naturalness, radial projection, and attainment of the minimum\\
\bottomrule
\end{tabularx}
\endgroup
\end{table}

The decisive difference arises when the radial function
$t\mapsto I(tu)$ may have several critical points or none. The mountain pass
method may still work if the global barrier and the
Palais--Smale condition remain in place. Nehari reduction, on the other hand, requires hypotheses
controlling the intersection of each ray with $\mathcal N$. When this
intersection is unique, Nehari's method can be viewed as a radial minimax; the mountain pass
method allows competition among all paths, not only rays.

\subsection{Exact equality of the levels}

For each direction $u$, one can construct a radial path ending at a
large multiple of $u$. If the class $\Gamma(0,e)$ has a fixed endpoint $e$,
this multiple must be connected to $e$. The region of negative energy allows
this connection to be made without increasing the maximum height of the path.

\begin{theorem}[Equality of the mountain pass and Nehari levels]
\label{teo:igualdad-niveles-montana-nehari-potencia}
Let $e\in X$ satisfy $I(e)<0$, and define
\[
 c_e:=\inf_{\gamma\in\Gamma(0,e)}
 \max_{t\in[0,1]}I(\gamma(t)).
\]
Then
\[
 c_e=m
 =\inf_{u\in\mathcal N}I(u)
 =\inf_{u\neq0}\max_{t>0}I(tu).
\]
In particular, the level does not depend on the choice of the endpoint with
negative energy.
\end{theorem}

\begin{proof}
We divide the proof into the two inequalities.

\emph{First inequality: $m\leq c_e$.}
Let $C_S$ be a constant for the embedding $X\hookrightarrow L^r(M)$. For
$v\in X$,
\[
 G(v)=\|v\|^p-\|v\|_{L^r(M)}^r
 \geq\|v\|^p-C_S^r\|v\|^r.
\]
Choose $\eta>0$ so small that $C_S^r\eta^{r-p}<1$. Then
\begin{equation}
 0<\|v\|\leq\eta
 \quad\Longrightarrow\quad G(v)>0.
\label{eq:G-positivo-cerca-origen}
\end{equation}
On the other hand, $I(e)<0$ implies
\[
 \frac1p\|e\|^p<\frac1r\|e\|_{L^r(M)}^r,
\]
and, since $r/p>1$,
\begin{equation}
 G(e)=\|e\|^p-\|e\|_{L^r(M)}^r<0.
\label{eq:G-negativo-extremo}
\end{equation}
Decrease $\eta$, if necessary, so that $\eta<\|e\|$.

Let $\gamma\in\Gamma(0,e)$. The set
\[
 A_\gamma:=\{t\in[0,1]\mid\|\gamma(t)\|=\eta\}
\]
is nonempty and compact. Let $t_0:=\displaystyle\max A_\gamma$. For $t>t_0$, we have
$\|\gamma(t)\|>\eta$: if the norm became less than or equal to $\eta$ again,
the intermediate value theorem would yield another point of $A_\gamma$ after
$t_0$. By \eqref{eq:G-positivo-cerca-origen},
$G(\gamma(t_0))>0$, whereas
\eqref{eq:G-negativo-extremo} gives $G(\gamma(1))<0$. Thus there exists
$t_*\in(t_0,1)$ such that $G(\gamma(t_*))=0$. Since
$\|\gamma(t_*)\|>\eta$, the point $\gamma(t_*)$ is nonzero and belongs to
$\mathcal N$. Consequently,
\[
 \max_{t\in[0,1]}I(\gamma(t))
 \geq I(\gamma(t_*))\geq m.
\]
Taking the infimum over $\gamma$ gives $c_e\geq m$.

\emph{Second inequality: $c_e\leq m$.}
Fix $u\in\mathcal N$ and set
\[
 a:=\frac{u}{\|u\|},
 \qquad
 b:=\frac{e}{\|e\|}.
\]
The unit sphere of $X$ is path connected. Indeed, if
$b\neq-a$, the path
\[
 s\longmapsto\frac{(1-s)a+sb}{\|(1-s)a+sb\|}
\]
joins $a$ to $b$; if $b=-a$, choose a unit vector
$z\notin\operatorname{span}\{a\}$ and concatenate the normalized paths from
$a$ to $z$ and from $z$ to $b$. Denote
one of these paths on the unit sphere by $h\colon[0,1]\to X$.

Such a vector $z$ exists because $X=W_0^{1,p}(M)$ is infinite-dimensional: a
sequence of nonzero smooth functions supported in disjoint balls
in the interior of $M$ is linearly independent. In particular,
$\dim X\geq2$.

The function $s\mapsto\|h(s)\|_{L^r(M)}^r$ is continuous and strictly
positive. By compactness of $[0,1]$, there exists
\[
 q:=\min_{s\in[0,1]}\|h(s)\|_{L^r(M)}^r>0.
\]
Choose $R>\displaystyle\max\{\|u\|,\|e\|\}$ so large that
\[
 \frac{R^p}{p}-\frac{qR^r}{r}<0.
\]
Then
\begin{equation}
 I(Rh(s))\leq\frac{R^p}{p}-\frac{qR^r}{r}<0
 \qquad(0\leq s\leq1).
\label{eq:conexion-negativa-esfera-grande}
\end{equation}

We construct $\gamma_u\in\Gamma(0,e)$ in three segments:
\[
 \gamma_u(t):=
 \begin{cases}
 3tR a,&0\leq t\leq\frac13,\\[1mm]
 R h(3t-1),&\frac13\leq t\leq\frac23,\\[1mm]
 \bigl((3-3t)R+(3t-2)\|e\|\bigr)b,
 &\frac23\leq t\leq1.
 \end{cases}
\]
The values agree at the joining points, because $h(0)=a$ and $h(1)=b$.
On the first segment, the radial function attains its unique maximum at
$\|u\|a=u$; thus the energy is at most $I(u)$. The second segment
has negative energy by
\eqref{eq:conexion-negativa-esfera-grande}. Finally, for
$s\geq\|e\|$,
\[
 I(sb)=\frac{s^p}{p}-\frac{s^r}{r}\|b\|_{L^r(M)}^r<0,
\]
since this inequality holds at $s=\|e\|$ and its equivalent form
$s^{r-p}>r/(p\|b\|_{L^r(M)}^r)$ remains valid as $s$ increases.
Thus the third segment also has negative energy, and
\[
 \max_{t\in[0,1]}I(\gamma_u(t))=I(u).
\]
It follows that $c_e\leq I(u)$ for every $u\in\mathcal N$. Taking the infimum gives
$c_e\leq m$. Equality with the radial minimax was already proved in
\eqref{eq:nivel-nehari-minimax-radial-potencia}.
\end{proof}

By Theorem~\ref{teo:nehari-potencia-minimizador-completo}, the common level
$c_e=m$ is attained at a ground state. By
Corollary~\ref{cor:paso-montana-forma-usual} and the (PS) condition established in
Proposition~\ref{prop:PS-potencia-p-paso-montana}, it is also a
critical level obtained by the mountain pass method. The two constructions need not
select the same function when there are several critical points, but
they do select the same minimum energy in this homogeneous model.

\section{A quasilinear equation on a vector bundle}
\label{sec:problema-quasilineal-haz-nehari}

We now apply the preceding theory to the quasilinear problem studied in \cite{VelazquezSandovalCarlosQuasilinear}. The characterization of the Dirichlet space by the trace allows us to interpret the homogeneous boundary condition, while the Sobolev and Rellich--Kondrashov embeddings provide the integrability and compactness required by the energy functional. We shall establish its Fréchet differentiability, the geometry of its radial functions, and the existence of a ground state by the Nehari method of Szulkin and Weth.

Let $(M,\mathbf{g})$ be a nonempty compact Riemannian manifold of dimension $n$, with smooth boundary, such that each connected component has nonempty boundary. Let $\mathbf{E}\to M$ be a smooth real vector bundle of positive finite rank, equipped with a bundle metric $\mathbf{h}_{\mathbf{E}}$ and a connection $\nabla^{\mathbf{E}}$ compatible with $\mathbf{h}_{\mathbf{E}}$. The metric $\mathbf{g}$ and the bundle metric $\mathbf{h}_{\mathbf{E}}$ induce the bundle metric $\mathbf{g}\otimes \mathbf{h}_{\mathbf{E}}$ on $T^*M\otimes \mathbf{E}$. Fix
\[
2\leq p\leq q<n,
\qquad
q^*:=\frac{nq}{n-q},
\]
and consider the problem
\begin{equation}\label{eq:problema-quasilineal-haz}
\begin{cases}
(\nabla^{\mathbf{E}})^*\!\left(a(|\nabla^{\mathbf{E}} \mathbf{u}|_{\mathbf{g}\otimes \mathbf{h}_{\mathbf{E}}}^{p})|\nabla^{\mathbf{E}} \mathbf{u}|_{\mathbf{g}\otimes \mathbf{h}_{\mathbf{E}}}^{p-2}\nabla^{\mathbf{E}} \mathbf{u}\right)=f(\mathbf{u})&\text{in }M,\\
\mathbf{u}=0&\text{on }\partial M.
\end{cases}
\end{equation}

Define $A(t):=\displaystyle\int_0^t a(s)\,ds$. We assume that:
\begin{enumerate}[label=\textup{(a\arabic*)},ref=a\arabic*]
\item\label{hip:quasilineal-a1} $a\colon [0,+\infty)\longrightarrow[0,+\infty)$ is continuous and there exist $k_1,k_2,k_3,k_4\geq0$, with $k_1>0$ and $k_2>0$ when $p<q$, such that
\[
k_1t^p+k_2t^q
\leq
a(t^p)t^p
\leq
k_3t^p+k_4t^q
\]
for every $t>0$;
\item\label{hip:quasilineal-a2} the function $\displaystyle t\longmapsto\displaystyle\frac{a(t^p)}{t^{q-p}}$ is nonincreasing on $(0,+\infty)$;
\item\label{hip:quasilineal-a3} the function $t\longmapsto A(t^p)$ is convex on $[0,+\infty)$.
\end{enumerate}

Let $\mathcal F\colon \mathbf{E}\longrightarrow\mathbb R$ be a continuous function whose restriction $\mathcal F_x\colon \mathbf{E}_x\longrightarrow\mathbb R$ to each fiber is of class $C^1$. We denote by $f_x\colon \mathbf{E}_x\longrightarrow \mathbf{E}_x$ its gradient with respect to $\mathbf{h}_{\mathbf{E}}$, that is,
\[
D\mathcal F_x(\xi)(\zeta)
=
\langle f_x(\xi),\zeta\rangle_{\mathbf{h}_{\mathbf{E}}}.
\]
We also assume that the bundle map $f\colon \mathbf{E}\longrightarrow \mathbf{E}$, defined by $f(\xi):=f_{\pi_{\mathbf{E}}(\xi)}(\xi)$, is continuous and that:
\begin{enumerate}[label=\textup{(f\arabic*)},ref=f\arabic*]
\item\label{hip:quasilineal-f1} $\displaystyle\lim_{|\xi|_{\mathbf{h}_{\mathbf{E}}}\to0}\displaystyle\frac{|f_x(\xi)|_{\mathbf{h}_{\mathbf{E}}}}{|\xi|_{\mathbf{h}_{\mathbf{E}}}^{p-1}}=0$ uniformly with respect to $x\in M$; that is, for each $\varepsilon>0$ there exists $\delta>0$, independent of $x$, such that
\[
|f_x(\xi)|_{\mathbf{h}_{\mathbf{E}}}
\leq
\varepsilon|\xi|_{\mathbf{h}_{\mathbf{E}}}^{p-1}
\]
if $x\in M$ and $0<|\xi|_{\mathbf{h}_{\mathbf{E}}}<\delta$;
\item\label{hip:quasilineal-f2} there exist $C_f>0$ and $r\in(q,q^*)$ such that
\[
|f_x(\xi)|_{\mathbf{h}_{\mathbf{E}}}
\leq
C_f\bigl(1+|\xi|_{\mathbf{h}_{\mathbf{E}}}^{r-1}\bigr)
\]
for every $x\in M$ and every $\xi\in \mathbf{E}_x$;
\item\label{hip:quasilineal-f3} there exists $\theta\in(q,r]$ such that
\[
0<\theta\mathcal F_x(\xi)
\leq
\langle f_x(\xi),\xi\rangle_{\mathbf{h}_{\mathbf{E}}}
\]
for every $x\in M$ and every $\xi\in \mathbf{E}_x\setminus\{0_x\}$;
\item\label{hip:quasilineal-f4} for each $x\in M$ and each $\xi\in \mathbf{E}_x\setminus\{0_x\}$, the function
\[
t\longmapsto
\frac{\langle f_x(t\xi),t\xi\rangle_{\mathbf{h}_{\mathbf{E}}}}{t^q}
\]
is strictly increasing on $(0,+\infty)$.
\end{enumerate}

The preceding hypotheses also determine the value of the potential on the zero section. Fix $x\in M$ and $\xi\in \mathbf{E}_x\setminus\{0_x\}$. By \textup{(\ref{hip:quasilineal-f3})}, for every $t>0$,
\[
0<\theta\mathcal F_x(t\xi)
\leq
\langle f_x(t\xi),t\xi\rangle_{\mathbf{h}_{\mathbf{E}}}.
\]
Let $\varepsilon>0$. By \textup{(\ref{hip:quasilineal-f1})}, there exists $\delta>0$ such that, if $0<t|\xi|_{\mathbf{h}_{\mathbf{E}}}<\delta$, then
\[
|f_x(t\xi)|_{\mathbf{h}_{\mathbf{E}}}
\leq
\varepsilon t^{p-1}|\xi|_{\mathbf{h}_{\mathbf{E}}}^{p-1}.
\]
Consequently,
\[
0<\mathcal F_x(t\xi)
\leq
\frac{\varepsilon}{\theta}t^p|\xi|_{\mathbf{h}_{\mathbf{E}}}^p
\]
for every sufficiently small $t>0$. Thus $\mathcal F_x(t\xi)\to0$ as $t\to0^+$. Since $\mathcal F_x$ is continuous on the fiber, we conclude that $\mathcal F_x(0_x)=0$. Theorem~\ref{teo:b5-fundamental-calculo-riemann-banach}, applied to $t\longmapsto\mathcal F_x(t\xi)$, then gives
\begin{equation}
\label{eq:representacion-potencial-fibrado}
\mathcal F_x(\xi)
=
\int_0^1\langle f_x(t\xi),\xi\rangle_{\mathbf{h}_{\mathbf{E}}}\,dt.
\end{equation}
Global continuity of $\mathcal F\colon \mathbf{E}\longrightarrow\mathbb R$ is part of the hypotheses. This condition does not follow merely from continuity of each restriction $\mathcal F_x$ and will be used, in particular, when taking minima and maxima over the compact unit sphere and ball bundles.

Take $X:=W_0^{1,q}(M,\mathbf{E})$. This space contains a nonzero section: simply multiply an element of a local frame of the bundle by a nonzero smooth function compactly supported in the interior of a trivializing chart. By the Poincaré inequality in Theorem~\ref{teo:poincare-W0-haces}, there exists a constant $C_P>0$, depending only on $M$, $\mathbf{g}$, $\mathbf{E}$, $\mathbf{h}_{\mathbf{E}}$, $\nabla^{\mathbf{E}}$, and $q$, such that, if
\[
\|\mathbf{u}\|_X
:=
\left(\int_M|\nabla_w^{\mathbf{E}}\mathbf{u}|_{\mathbf{g}\otimes \mathbf{h}_{\mathbf{E}}}^q\,d\lambda_{\mathbf{g}}\right)^{\frac{1}{q}},
\]
then
\[
\|\mathbf{u}\|_X
\leq
\|\mathbf{u}\|_{W^{1,q}(M,\mathbf{E})}
\leq
(1+C_P^q)^{\frac{1}{q}}\|\mathbf{u}\|_X
\]
for every $\mathbf{u}\in X$. We use $\|\cdot\|_X$ as the norm on $X$.
Henceforth, we write $\nabla^{\mathbf{E}} \mathbf{u}$ for the weak covariant derivative of $\mathbf{u}\in X$. Since $p\leq q$ and the volume is finite, Hölder's inequality gives $W_0^{1,q}(M,\mathbf{E})\hookrightarrow W_0^{1,p}(M,\mathbf{E})$: a sequence of compactly supported smooth sections converging in $W^{1,q}(M,\mathbf E)$ also converges in $W^{1,p}(M,\mathbf E)$. Thus Poincaré's inequality applies to both exponents $p$ and $q$.

\begin{definition}[Weak solution of the quasilinear problem]
\label{def:solucion-debil-problema-quasilineal-haz}
Let $M$ be the compact Riemannian manifold with nonempty boundary fixed in this
section, and let $\mathbf{E}\to M$ be the vector bundle specified there.
A section $\mathbf{u}\in X$ is a weak solution of \eqref{eq:problema-quasilineal-haz} if
\[
\int_M
a(|\nabla^{\mathbf{E}} \mathbf{u}|_{\mathbf{g}\otimes \mathbf{h}_{\mathbf{E}}}^{p})
|\nabla^{\mathbf{E}} \mathbf{u}|_{\mathbf{g}\otimes \mathbf{h}_{\mathbf{E}}}^{p-2}
\langle\nabla^{\mathbf{E}} \mathbf{u},\nabla^{\mathbf{E}}\boldsymbol{\phi}\rangle_{\mathbf{g}\otimes \mathbf{h}_{\mathbf{E}}}\,d\lambda_{\mathbf{g}}
=
\int_M\langle f_x(\mathbf{u}),\boldsymbol{\phi}\rangle_{\mathbf{h}_{\mathbf{E}}}\,d\lambda_{\mathbf{g}}
\]
for every $\boldsymbol{\phi}\in X$.
\end{definition}

\subsection{Preliminary estimates}

The monotonicity in \textup{(\ref{hip:quasilineal-a2})} relates the primitive $A$ to the term appearing in the Nehari identity.

\begin{lemma}
\label{lem:desigualdad-A-a-quasilineal}
For every $\rho\geq0$, we have
\[
A(\rho^p)
\geq
\frac{p}{q}a(\rho^p)\rho^p.
\]
\end{lemma}

\begin{proof}
The assertion is immediate for $\rho=0$. Let $\rho>0$. If $0<\tau\leq\rho$, hypothesis \textup{(\ref{hip:quasilineal-a2})} implies $\displaystyle\frac{a(\tau^p)}{\tau^{q-p}}
\geq
\displaystyle\frac{a(\rho^p)}{\rho^{q-p}}$, and hence $a(\tau^p)
\geq
a(\rho^p)\left(\displaystyle\frac{\tau}{\rho}\right)^{q-p}$.

Using the change of variable $s=\tau^p$, we obtain
\[
A(\rho^p)
=
p\int_0^\rho a(\tau^p)\tau^{p-1}\,d\tau
\geq
p\,a(\rho^p)\rho^{p-q}\int_0^\rho\tau^{q-1}\,d\tau
=
\frac{p}{q}a(\rho^p)\rho^p.
\]
\end{proof}

\begin{lemma}[Fiberwise growth estimates]
\label{lem:estimaciones-crecimiento-fibra-quasilineal}
Let $M$ be the compact Riemannian manifold with nonempty boundary fixed in this
section, and let $\mathbf{E}\to M$ be the vector bundle specified there.
For each $\varepsilon>0$ there exists $C_\varepsilon>0$ such that
\begin{equation}\label{eq:estimacion-f-quasilineal}
|f_x(\xi)|_{\mathbf{h}_{\mathbf{E}}}
\leq
\varepsilon|\xi|_{\mathbf{h}_{\mathbf{E}}}^{p-1}
+
C_\varepsilon|\xi|_{\mathbf{h}_{\mathbf{E}}}^{r-1}
\end{equation}
for every $x\in M$ and every $\xi\in \mathbf{E}_x$. Moreover,
\begin{equation}\label{eq:estimacion-F-quasilineal}
|\mathcal F_x(\xi)|
\leq
\frac{\varepsilon}{p}|\xi|_{\mathbf{h}_{\mathbf{E}}}^{p}
+
\frac{C_\varepsilon}{r}|\xi|_{\mathbf{h}_{\mathbf{E}}}^{r}.
\end{equation}
\end{lemma}

\begin{proof}
Let $\varepsilon>0$. By \textup{(\ref{hip:quasilineal-f1})}, there exists $\delta\in(0,1)$ such that
\[
|f_x(\xi)|_{\mathbf{h}_{\mathbf{E}}}
\leq
\varepsilon|\xi|_{\mathbf{h}_{\mathbf{E}}}^{p-1}
\]
if $x\in M$ and $0<|\xi|_{\mathbf{h}_{\mathbf{E}}}\leq\delta$. Letting $\xi\to0_x$ and using continuity of $f_x$ gives $f_x(0_x)=0$, so the same inequality holds for $\xi=0_x$.

If $|\xi|_{\mathbf{h}_{\mathbf{E}}}\geq\delta$, hypothesis \textup{(\ref{hip:quasilineal-f2})} gives
\[
|f_x(\xi)|_{\mathbf{h}_{\mathbf{E}}}
\leq
C_f\bigl(1+|\xi|_{\mathbf{h}_{\mathbf{E}}}^{r-1}\bigr)
\leq
C_f\bigl(1+\delta^{1-r}\bigr)|\xi|_{\mathbf{h}_{\mathbf{E}}}^{r-1}.
\]
Taking $C_\varepsilon\geq C_f(1+\delta^{1-r})$ yields \eqref{eq:estimacion-f-quasilineal} in both regions.

Since $\mathcal F_x(0_x)=0$ and $f_x$ is the fiberwise gradient of $\mathcal F_x$, Theorem~\ref{teo:b5-fundamental-calculo-riemann-banach}, applied to $t\longmapsto\mathcal F_x(t\xi)$, gives
\[
\mathcal F_x(\xi)
=
\int_0^1\langle f_x(t\xi),\xi\rangle_{\mathbf{h}_{\mathbf{E}}}\,dt.
\]
By \eqref{eq:estimacion-f-quasilineal},
\[
|\mathcal F_x(\xi)|\leq
\int_0^1|f_x(t\xi)|_{\mathbf{h}_{\mathbf{E}}}|\xi|_{\mathbf{h}_{\mathbf{E}}}\,dt\leq
\int_0^1
\left(
\varepsilon t^{p-1}|\xi|_{\mathbf{h}_{\mathbf{E}}}^{p-1}
+
C_\varepsilon t^{r-1}|\xi|_{\mathbf{h}_{\mathbf{E}}}^{r-1}
\right)|\xi|_{\mathbf{h}_{\mathbf{E}}}\,dt=
\frac{\varepsilon}{p}|\xi|_{\mathbf{h}_{\mathbf{E}}}^{p}
+
\frac{C_\varepsilon}{r}|\xi|_{\mathbf{h}_{\mathbf{E}}}^{r}.
\]
\end{proof}

\begin{lemma}[Lower bound for the potential]
\label{lem:cota-inferior-potencial-quasilineal}
Let $M$ be the compact Riemannian manifold with nonempty boundary fixed in this
section, and let $\mathbf{E}\to M$ be the vector bundle specified there.
There exist constants $D>0$ and $D'\geq0$ such that
\[
\mathcal F_x(\xi)
\geq
D|\xi|_{\mathbf{h}_{\mathbf{E}}}^{\theta}-D'
\]
for every $x\in M$ and every $\xi\in \mathbf{E}_x$.
\end{lemma}

\begin{proof}
Fix $x\in M$ and $\eta\in \mathbf{E}_x\setminus\{0_x\}$. For $t>0$, define
\[
H_{x,\eta}(t):=\frac{\mathcal F_x(t\eta)}{t^\theta}.
\]
By \textup{(\ref{hip:quasilineal-f3})}, $\mathcal F_x(t\eta)>0$. Moreover,
\[
\frac{d}{dt}\log H_{x,\eta}(t)
=
\frac{\langle f_x(t\eta),t\eta\rangle_{\mathbf{h}_{\mathbf{E}}}-\theta\mathcal F_x(t\eta)}
{t\mathcal F_x(t\eta)}
\geq0.
\]
Thus $H_{x,\eta}$ is nondecreasing.

The unit sphere bundle
\[
S(\mathbf{E}):=\{\eta\in \mathbf{E}\mid |\eta|_{\mathbf{h}_{\mathbf{E}}}=1\}
\]
is compact, because $M$ is compact and $\mathbf{E}$ has finite rank. The function $\eta\longmapsto\mathcal F_{\pi_{\mathbf{E}}(\eta)}(\eta)$ is continuous and, by \textup{(\ref{hip:quasilineal-f3})}, positive on $S(\mathbf{E})$. Consequently,
\[
D:=\min_{\eta\in S(\mathbf{E})}\mathcal F_{\pi_{\mathbf{E}}(\eta)}(\eta)>0.
\]
If $|\xi|_{\mathbf{h}_{\mathbf{E}}}\geq1$, write $\xi=\rho\eta$, with $\rho=|\xi|_{\mathbf{h}_{\mathbf{E}}}\geq1$ and $|\eta|_{\mathbf{h}_{\mathbf{E}}}=1$. Monotonicity of $H_{x,\eta}$ implies
\[
\mathcal F_x(\xi)
=
\mathcal F_x(\rho\eta)
\geq
\rho^\theta\mathcal F_x(\eta)
\geq
D|\xi|_{\mathbf{h}_{\mathbf{E}}}^{\theta}.
\]

The closed unit ball bundle $\{\xi\in \mathbf{E}\mid |\xi|_{\mathbf{h}_{\mathbf{E}}}\leq1\}$ is compact. Thus the continuous function
\[
\xi\longmapsto D|\xi|_{\mathbf{h}_{\mathbf{E}}}^{\theta}-\mathcal F_{\pi_{\mathbf{E}}(\xi)}(\xi)
\]
attains a maximum $D'\geq0$ on this set. Hence $\mathcal F_x(\xi)
\geq
D|\xi|_{\mathbf{h}_{\mathbf{E}}}^{\theta}-D'$ also when $|\xi|_{\mathbf{h}_{\mathbf{E}}}\leq1$.
\end{proof}

\subsection{The energy functional}

Define
\begin{equation}\label{eq:funcional-energia-quasilineal}
I(\mathbf{u})
:=
\frac{1}{p}\int_M A(|\nabla^{\mathbf{E}} \mathbf{u}|_{\mathbf{g}\otimes \mathbf{h}_{\mathbf{E}}}^p)\,d\lambda_{\mathbf{g}}
-
\int_M\mathcal F_x(\mathbf{u}(x))\,d\lambda_{\mathbf{g}},
\qquad \mathbf{u}\in X.
\end{equation}

\begin{proposition}[Differentiability and weak solutions]
\label{prop:energia-quasilineal-C1}
Let $M$ be the compact Riemannian manifold with nonempty boundary fixed in this
section, and let $\mathbf{E}\to M$ be the vector bundle specified there.
The functional $I\colon X\longrightarrow\mathbb R$ is well defined and belongs to $C^1(X,\mathbb R)$. For all $\mathbf{u},\mathbf{v}\in X$,
\begin{equation}\label{eq:derivada-energia-quasilineal}
DI(\mathbf{u})(\mathbf{v})=
\int_M a(|\nabla^{\mathbf{E}} \mathbf{u}|_{\mathbf{g}\otimes \mathbf{h}_{\mathbf{E}}}^{p})
|\nabla^{\mathbf{E}} \mathbf{u}|_{\mathbf{g}\otimes \mathbf{h}_{\mathbf{E}}}^{p-2}
\langle\nabla^{\mathbf{E}} \mathbf{u},\nabla^{\mathbf{E}} \mathbf{v}\rangle_{\mathbf{g}\otimes \mathbf{h}_{\mathbf{E}}}\,d\lambda_{\mathbf{g}}-
\int_M\langle f_x(\mathbf{u}(x)),\mathbf{v}(x)\rangle_{\mathbf{h}_{\mathbf{E}}}\,d\lambda_{\mathbf{g}}.
\end{equation}
Consequently, the critical points of $I$ are exactly the weak solutions of \eqref{eq:problema-quasilineal-haz}.
\end{proposition}

\begin{proof}
Write $I=J_1-J_2$, where
\[
J_1(\mathbf{u}):=\frac{1}{p}\int_M A(|\nabla^{\mathbf{E}} \mathbf{u}|_{\mathbf{g}\otimes \mathbf{h}_{\mathbf{E}}}^p)\,d\lambda_{\mathbf{g}},
\qquad
J_2(\mathbf{u}):=\int_M\mathcal F_x(\mathbf{u}(x))\,d\lambda_{\mathbf{g}}.
\]

We first check that both terms are well defined. By \textup{(\ref{hip:quasilineal-a1})}, for each $\rho>0$ we have $a(\rho^p)\leq k_3+k_4\rho^{q-p}$. Thus
\[
A(\rho^p)
=
p\int_0^\rho a(\tau^p)\tau^{p-1}\,d\tau
\leq
k_3\rho^p+\frac{pk_4}{q}\rho^q.
\]
Since $\nabla^{\mathbf{E}}\mathbf{u}\in L^q(M,T^*M\otimes \mathbf{E})$, $p\leq q$, and $\lambda_{\mathbf{g}}(M)<\infty$, we obtain $A(|\nabla^{\mathbf{E}}\mathbf{u}|_{\mathbf{g}\otimes \mathbf{h}_{\mathbf{E}}}^p)\in L^1(M)$.

By Lemma~\ref{lem:estimaciones-crecimiento-fibra-quasilineal}, there exists $C>0$ such that $|\mathcal F_x(\xi)|
\leq
C\bigl(|\xi|_{\mathbf{h}_{\mathbf{E}}}^p+|\xi|_{\mathbf{h}_{\mathbf{E}}}^r\bigr)$.

The continuous embeddings $X\hookrightarrow L^p(M,\mathbf{E})$ and $X\hookrightarrow L^r(M,\mathbf{E})$ show that $\mathcal F_x(\mathbf{u}(x))\in L^1(M)$. Thus $I$ is well defined.

We now prove that $J_1\in C^1(X,\mathbb R)$. Let $\mathbf{V}:=T^*M\otimes \mathbf{E}$, equipped with the bundle metric $\mathbf{g}\otimes \mathbf{h}_{\mathbf{E}}$, and define the bundle maps
\[
\Psi_x\colon \mathbf{V}_x\longrightarrow \mathbf{V}_x,
\qquad
\Psi_x(\eta)
:=
a(|\eta|_{\mathbf{g}\otimes \mathbf{h}_{\mathbf{E}}}^p)|\eta|_{\mathbf{g}\otimes \mathbf{h}_{\mathbf{E}}}^{p-2}\eta,
\qquad
\Psi_x(0_x):=0_x,
\]
and $\Phi_x\colon \mathbf{V}_x\longrightarrow\mathbb R$, $\Phi_x(\eta):=\displaystyle\frac{1}{p}A(|\eta|_{\mathbf{g}\otimes \mathbf{h}_{\mathbf{E}}}^p)$.

The map $\Psi\colon \mathbf{V}\longrightarrow \mathbf{V}$ is continuous and
\[
D\Phi_x(\eta)(\zeta)
=
\langle\Psi_x(\eta),\zeta\rangle_{\mathbf{g}\otimes \mathbf{h}_{\mathbf{E}}}.
\]
Moreover, by \textup{(\ref{hip:quasilineal-a1})},
\[
|\Psi_x(\eta)|_{\mathbf{g}\otimes \mathbf{h}_{\mathbf{E}}}
\leq
k_3|\eta|_{\mathbf{g}\otimes \mathbf{h}_{\mathbf{E}}}^{p-1}
+k_4|\eta|_{\mathbf{g}\otimes \mathbf{h}_{\mathbf{E}}}^{q-1}
\leq
C\bigl(1+|\eta|_{\mathbf{g}\otimes \mathbf{h}_{\mathbf{E}}}^{q-1}\bigr).
\]
Corollary~\ref{cor:superposicion-exponentes-conjugados-haces}, applied with $m=q$, implies that $\eta\longmapsto\Psi(\eta)$ defines a continuous operator from $L^q(M,\mathbf{V})$ to $L^{q'}(M,\mathbf{V})$, where $q':=\displaystyle\frac{q}{q-1}$.

For $\mathbf{u}\in X$, define $L_{\mathbf{u}}\in X'$ by
\[
L_{\mathbf{u}}(\mathbf{v})
:=
\int_M\langle\Psi_x(\nabla^{\mathbf{E}}\mathbf{u}),\nabla^{\mathbf{E}}\mathbf{v}\rangle_{\mathbf{g}\otimes \mathbf{h}_{\mathbf{E}}}\,d\lambda_{\mathbf{g}}.
\]
Hölder's inequality (Proposition~\ref{desigualdad de holder}) gives
\[
|L_{\mathbf{u}}(\mathbf{v})|
\leq
\|\Psi(\nabla^{\mathbf{E}}\mathbf{u})\|_{L^{q'}(M,\mathbf{V})}
\|\nabla^{\mathbf{E}}\mathbf{v}\|_{L^q(M,\mathbf{V})}
=
\|\Psi(\nabla^{\mathbf{E}}\mathbf{u})\|_{L^{q'}(M,\mathbf{V})}\|\mathbf{v}\|_X,
\]
so $L_{\mathbf{u}}$ is linear and continuous.

Let $\mathbf{h}\in X$. For almost every $x\in M$, Theorem~\ref{teo:b5-fundamental-calculo-riemann-banach}, applied to
\[
s\longmapsto\Phi_x(\nabla^{\mathbf{E}}\mathbf{u}(x)+s\nabla^{\mathbf{E}}\mathbf{h}(x))
\]
gives
\[
\Phi_x(\nabla^{\mathbf{E}}\mathbf{u}+\nabla^{\mathbf{E}}\mathbf{h})-\Phi_x(\nabla^{\mathbf{E}}\mathbf{u})
=
\int_0^1
\langle\Psi_x(\nabla^{\mathbf{E}}\mathbf{u}+s\nabla^{\mathbf{E}}\mathbf{h}),\nabla^{\mathbf{E}}\mathbf{h}\rangle_{\mathbf{g}\otimes \mathbf{h}_{\mathbf{E}}}\,ds.
\]
Hölder's inequality (Proposition~\ref{desigualdad de holder}) and the preceding bounds show that the integrand belongs to $L^1(M\times[0,1])$. We may apply Fubini's theorem for Bochner integrals~\ref{teo:b5-fubini-bochner}, in its scalar case, to $(M,\lambda_{\mathbf{g}})\times([0,1],ds)$ and obtain
\[
\begin{split}
J_1(\mathbf{u}+\mathbf{h})-J_1(\mathbf{u})-L_{\mathbf{u}}(\mathbf{h})
&=
\int_0^1\int_M
\langle
\Psi_x(\nabla^{\mathbf{E}}\mathbf{u}+s\nabla^{\mathbf{E}}\mathbf{h})-\Psi_x(\nabla^{\mathbf{E}}\mathbf{u}),
\nabla^{\mathbf{E}}\mathbf{h}
\rangle_{\mathbf{g}\otimes \mathbf{h}_{\mathbf{E}}}\,d\lambda_{\mathbf{g}}\,ds.
\end{split}
\]
Consequently,
\[
\begin{split}
|J_1(\mathbf{u}+\mathbf{h})-J_1(\mathbf{u})-L_{\mathbf{u}}(\mathbf{h})|
&\leq
\|\mathbf{h}\|_X
\int_0^1
\|\Psi(\nabla^{\mathbf{E}}\mathbf{u}+s\nabla^{\mathbf{E}}\mathbf{h})-\Psi(\nabla^{\mathbf{E}}\mathbf{u})\|_{L^{q'}(M,\mathbf{V})}\,ds.
\end{split}
\]
Let $(\mathbf{h}_j)_{j\in\mathbb N}\subseteq X$ satisfy $\mathbf{h}_j\to0$ in $X$. For each $s\in[0,1]$,
\[
\nabla^{\mathbf{E}}\mathbf{u}+s\nabla^{\mathbf{E}}\mathbf{h}_j
\longrightarrow
\nabla^{\mathbf{E}}\mathbf{u}
\quad\text{in }L^q(M,\mathbf{V}),
\]
and continuity of the superposition operator induced by $\Psi$ implies
\[
\|\Psi(\nabla^{\mathbf{E}}\mathbf{u}+s\nabla^{\mathbf{E}}\mathbf{h}_j)-\Psi(\nabla^{\mathbf{E}}\mathbf{u})\|_{L^{q'}(M,\mathbf{V})}
\longrightarrow0.
\]
Since $(\mathbf{h}_j)$ is bounded in $X$, there exists $R>0$ such that $\|\nabla^{\mathbf{E}}\mathbf{u}+s\nabla^{\mathbf{E}}\mathbf{h}_j\|_{L^q(M,\mathbf{V})}\leq R$ for every $j$ and every $s\in[0,1]$. The growth bound for $\Psi$ provides a constant $C_R>0$ such that
\[
\|\Psi(\nabla^{\mathbf{E}}\mathbf{u}+s\nabla^{\mathbf{E}}\mathbf{h}_j)\|_{L^{q'}(M,\mathbf{V})}
+
\|\Psi(\nabla^{\mathbf{E}}\mathbf{u})\|_{L^{q'}(M,\mathbf{V})}
\leq C_R
\]
for every $j$ and every $s\in[0,1]$. Theorem~\ref{convergencia dominada}, applied in the variable $s$, implies
\[
\int_0^1
\|\Psi(\nabla^{\mathbf{E}}\mathbf{u}+s\nabla^{\mathbf{E}}\mathbf{h}_j)-\Psi(\nabla^{\mathbf{E}}\mathbf{u})\|_{L^{q'}(M,\mathbf{V})}\,ds
\longrightarrow0.
\]
Thus
\[
\lim_{\mathbf{h}\to0}
\frac{|J_1(\mathbf{u}+\mathbf{h})-J_1(\mathbf{u})-L_{\mathbf{u}}(\mathbf{h})|}{\|\mathbf{h}\|_X}
=0.
\]
Hence $J_1$ is Fréchet differentiable and $DJ_1(\mathbf{u})=L_{\mathbf{u}}$.

If $\mathbf{u}_j\to \mathbf{u}$ in $X$, then $\nabla^{\mathbf{E}}\mathbf{u}_j\to\nabla^{\mathbf{E}}\mathbf{u}$ in $L^q(M,\mathbf{V})$ and
\[
\Psi(\nabla^{\mathbf{E}}\mathbf{u}_j)\to\Psi(\nabla^{\mathbf{E}}\mathbf{u})
\quad\text{in }L^{q'}(M,\mathbf{V}).
\]
Consequently,
\[
\|DJ_1(\mathbf{u}_j)-DJ_1(\mathbf{u})\|_{X'}
\leq
\|\Psi(\nabla^{\mathbf{E}}\mathbf{u}_j)-\Psi(\nabla^{\mathbf{E}}\mathbf{u})\|_{L^{q'}(M,\mathbf{V})}
\longrightarrow0.
\]
This proves that $J_1\in C^1(X,\mathbb R)$.

We now prove that $J_2\in C^1(X,\mathbb R)$. By \eqref{eq:estimacion-f-quasilineal} and $p<r$, there exists $C>0$ such that
\[
|f_x(\xi)|_{\mathbf{h}_{\mathbf{E}}}
\leq
C\bigl(1+|\xi|_{\mathbf{h}_{\mathbf{E}}}^{r-1}\bigr).
\]
Corollary~\ref{cor:superposicion-exponentes-conjugados-haces}, applied with $m=r$, implies that the operator $\mathbf{u}\longmapsto f(\mathbf{u})$ is continuous from $L^r(M,\mathbf{E})$ to $L^{r'}(M,\mathbf{E})$, where $r':=\displaystyle\frac{r}{r-1}$.

For $\mathbf{u}\in X$, define $K_{\mathbf{u}}\in X'$ by
\[
K_{\mathbf{u}}(\mathbf{v}):=\int_M\langle f_x(\mathbf{u}(x)),\mathbf{v}(x)\rangle_{\mathbf{h}_{\mathbf{E}}}\,d\lambda_{\mathbf{g}}.
\]
By Hölder's inequality (Proposition~\ref{desigualdad de holder}) and the embedding $X\hookrightarrow L^r(M,\mathbf{E})$, there exists a constant $C_S>0$ such that
\[
|K_{\mathbf{u}}(\mathbf{v})|
\leq
\|f(\mathbf{u})\|_{L^{r'}(M,\mathbf{E})}\|\mathbf{v}\|_{L^r(M,\mathbf{E})}
\leq
C_S\|f(\mathbf{u})\|_{L^{r'}(M,\mathbf{E})}\|\mathbf{v}\|_X.
\]

For $\mathbf{h}\in X$, Theorem~\ref{teo:b5-fundamental-calculo-riemann-banach}, applied in each fiber, gives
\[
\mathcal F_x(\mathbf{u}+\mathbf{h})-\mathcal F_x(\mathbf{u})
=
\int_0^1\langle f_x(\mathbf{u}+s\mathbf{h}),\mathbf{h}\rangle_{\mathbf{h}_{\mathbf{E}}}\,ds.
\]
The embeddings and the growth bound allow us to apply Fubini's theorem for Bochner integrals~\ref{teo:b5-fubini-bochner}, in its scalar case, to $(M,\lambda_{\mathbf{g}})\times([0,1],ds)$. Thus
\[
J_2(\mathbf{u}+\mathbf{h})-J_2(\mathbf{u})-K_{\mathbf{u}}(\mathbf{h})
=
\int_0^1\int_M
\langle f_x(\mathbf{u}+s\mathbf{h})-f_x(\mathbf{u}),\mathbf{h}\rangle_{\mathbf{h}_{\mathbf{E}}}\,d\lambda_{\mathbf{g}}\,ds,
\]
and
\[
|J_2(\mathbf{u}+\mathbf{h})-J_2(\mathbf{u})-K_{\mathbf{u}}(\mathbf{h})|
\leq
C_S\|\mathbf{h}\|_X
\int_0^1\|f(\mathbf{u}+s\mathbf{h})-f(\mathbf{u})\|_{L^{r'}(M,\mathbf{E})}\,ds.
\]
If $\mathbf{h}_j\to0$ in $X$, then $\mathbf{u}+s\mathbf{h}_j\to \mathbf{u}$ in $L^r(M,\mathbf{E})$ for each $s\in[0,1]$. Continuity of the superposition operator induced by $f$ implies
\[
\|f(\mathbf{u}+s\mathbf{h}_j)-f(\mathbf{u})\|_{L^{r'}(M,\mathbf{E})}\longrightarrow0.
\]
Since $(\mathbf{h}_j)$ is bounded in $X$ and the embedding $X\hookrightarrow L^r(M,\mathbf{E})$ is continuous, there exists $R>0$ such that
\[
\|\mathbf{u}+s\mathbf{h}_j\|_{L^r(M,\mathbf{E})}\leq R
\]
for every $j$ and every $s\in[0,1]$. The boundedness estimate in Remark~\ref{rem: nemytski-variedad-compacta-haz} provides $C_R>0$ such that
\[
\|f(\mathbf{u}+s\mathbf{h}_j)-f(\mathbf{u})\|_{L^{r'}(M,\mathbf{E})}
\leq
\|f(\mathbf{u}+s\mathbf{h}_j)\|_{L^{r'}(M,\mathbf{E})}+\|f(\mathbf{u})\|_{L^{r'}(M,\mathbf{E})}
\leq C_R
\]
for every $j$ and every $s\in[0,1]$. By Theorem~\ref{convergencia dominada}, applied in the variable $s$,
\[
\int_0^1\|f(\mathbf{u}+s\mathbf{h}_j)-f(\mathbf{u})\|_{L^{r'}(M,\mathbf{E})}\,ds
\longrightarrow0.
\]
It follows that $J_2$ is Fréchet differentiable and $DJ_2(\mathbf{u})=K_{\mathbf{u}}$.

Finally, if $\mathbf{u}_j\to \mathbf{u}$ in $X$, then $\mathbf{u}_j\to \mathbf{u}$ in $L^r(M,\mathbf{E})$ and $f(\mathbf{u}_j)\to f(\mathbf{u})$ in $L^{r'}(M,\mathbf{E})$. Hence
\[
\|DJ_2(\mathbf{u}_j)-DJ_2(\mathbf{u})\|_{X'}
\leq
C_S\|f(\mathbf{u}_j)-f(\mathbf{u})\|_{L^{r'}(M,\mathbf{E})}
\longrightarrow0.
\]
Therefore $J_2\in C^1(X,\mathbb R)$. Identity \eqref{eq:derivada-energia-quasilineal} follows from $I=J_1-J_2$. Comparing it with Definition~\ref{def:solucion-debil-problema-quasilineal-haz}, we conclude that $DI(\mathbf{u})=0$ if and only if $\mathbf{u}$ is a weak solution.
\end{proof}

\subsection{Radial geometry and the Nehari manifold}

Define
\[
\mathcal N
:=
\{\mathbf{u}\in X\setminus\{0\}\mid DI(\mathbf{u})(\mathbf{u})=0\}.
\]
Thus $\mathbf{u}\in\mathcal N$ if and only if $\mathbf{u}\neq0$ and
\begin{equation}\label{eq:identidad-nehari-quasilineal}
\int_M a(|\nabla^{\mathbf{E}}\mathbf{u}|_{\mathbf{g}\otimes \mathbf{h}_{\mathbf{E}}}^p)|\nabla^{\mathbf{E}}\mathbf{u}|_{\mathbf{g}\otimes \mathbf{h}_{\mathbf{E}}}^p\,d\lambda_{\mathbf{g}}
=
\int_M\langle f_x(\mathbf{u}),\mathbf{u}\rangle_{\mathbf{h}_{\mathbf{E}}}\,d\lambda_{\mathbf{g}}.
\end{equation}

\begin{proposition}[Projection onto the Nehari manifold]
\label{prop:proyeccion-nehari-quasilineal}
Let $M$ be the compact Riemannian manifold with nonempty boundary fixed in this
section, and let $\mathbf{E}\to M$ be the vector bundle specified there.
For each $\mathbf{u}\in X\setminus\{0\}$ there exists a unique $t_{\mathbf{u}}>0$ such that $t_{\mathbf{u}}\mathbf{u}\in\mathcal N$. Moreover, $t\longmapsto I(t\mathbf{u})$ is strictly increasing on $(0,t_{\mathbf{u}})$ and strictly decreasing on $(t_{\mathbf{u}},+\infty)$.
\end{proposition}

\begin{proof}
Fix $\mathbf{u}\in X\setminus\{0\}$ and define $\gamma_{\mathbf{u}}(t):=I(t\mathbf{u})$ for $t>0$. By the chain rule,
\[
\gamma_{\mathbf{u}}'(t)=DI(t\mathbf{u})(\mathbf{u})=\frac{1}{t}DI(t\mathbf{u})(t\mathbf{u}).
\]

We first show that $\gamma_{\mathbf{u}}(t)>0$ for sufficiently small $t>0$. By Lemma~\ref{lem:desigualdad-A-a-quasilineal} and \textup{(\ref{hip:quasilineal-a1})},
\[
\frac{1}{p}A(t^p|\nabla^{\mathbf{E}}\mathbf{u}|_{\mathbf{g}\otimes \mathbf{h}_{\mathbf{E}}}^p)
\geq
\frac{1}{q} a(t^p|\nabla^{\mathbf{E}}\mathbf{u}|_{\mathbf{g}\otimes \mathbf{h}_{\mathbf{E}}}^p)
t^p|\nabla^{\mathbf{E}}\mathbf{u}|_{\mathbf{g}\otimes \mathbf{h}_{\mathbf{E}}}^p
\geq
\frac{k_1}{q}t^p|\nabla^{\mathbf{E}}\mathbf{u}|_{\mathbf{g}\otimes \mathbf{h}_{\mathbf{E}}}^p+
\frac{k_2}{q}t^q|\nabla^{\mathbf{E}}\mathbf{u}|_{\mathbf{g}\otimes \mathbf{h}_{\mathbf{E}}}^q,
\]
where the norms of the derivatives are taken with respect to $\mathbf{g}\otimes \mathbf{h}_{\mathbf{E}}$. By \eqref{eq:estimacion-F-quasilineal},
\[
\mathcal F_x(t\mathbf{u})
\leq
\frac{\varepsilon}{p}t^p|\mathbf{u}|_{\mathbf{h}_{\mathbf{E}}}^p
+
\frac{C_\varepsilon}{r}t^r|\mathbf{u}|_{\mathbf{h}_{\mathbf{E}}}^r.
\]
Let $C_P>0$ satisfy $\|\mathbf{u}\|_{L^p(M,\mathbf{E})}\leq C_P\|\nabla^{\mathbf{E}}\mathbf{u}\|_{L^p(M,T^*M\otimes \mathbf{E})}$. Take $\varepsilon:=\displaystyle\frac{pk_1}{2qC_P^p}$.

Integrating the preceding inequalities gives
\[
\gamma_{\mathbf{u}}(t)
\geq
\frac{k_1}{2q}t^p\|\nabla^{\mathbf{E}}\mathbf{u}\|_{L^p(M,T^*M\otimes \mathbf{E})}^p
+
\frac{k_2}{q}t^q\|\nabla^{\mathbf{E}}\mathbf{u}\|_{L^q(M,T^*M\otimes \mathbf{E})}^q
-
\frac{C_\varepsilon}{r}t^r\|\mathbf{u}\|_{L^r(M,\mathbf{E})}^r.
\]
Since $\mathbf{u}\neq0$, Poincaré's inequality implies $\|\nabla^{\mathbf{E}}\mathbf{u}\|_{L^p(M,T^*M\otimes \mathbf{E})}>0$. Moreover, $r>q\geq p$. Thus there exists $t_0>0$ such that $\gamma_{\mathbf{u}}(t)>0$ for every $t\in(0,t_0)$.

On the other hand, the upper bound in \textup{(\ref{hip:quasilineal-a1})} implies $A(\rho^p)
\leq
k_3\rho^p+\displaystyle\frac{pk_4}{q}\rho^q$.

Lemma~\ref{lem:cota-inferior-potencial-quasilineal} gives
\[
\gamma_{\mathbf{u}}(t)
\leq
\frac{k_3}{p}t^p\|\nabla^{\mathbf{E}}\mathbf{u}\|_{L^p(M,T^*M\otimes \mathbf{E})}^p
+
\frac{k_4}{q}t^q\|\nabla^{\mathbf{E}}\mathbf{u}\|_{L^q(M,T^*M\otimes \mathbf{E})}^q
-Dt^\theta\|\mathbf{u}\|_{L^\theta(M,\mathbf{E})}^\theta
+D'\lambda_{\mathbf{g}}(M).
\]
Since $\theta>q\geq p$ and $\mathbf{u}\neq0$, we have $\gamma_{\mathbf{u}}(t)\to-\infty$ as $t\to+\infty$. We must also rule out a maximizing sequence escaping toward $0$. The upper bounds in \textup{(\ref{hip:quasilineal-a1})} and \eqref{eq:estimacion-F-quasilineal} give, for $0<t\leq1$, a constant $C_{\mathbf{u}}>0$ such that
\[
 |\gamma_{\mathbf{u}}(t)|
 \leq C_{\mathbf{u}}(t^p+t^q+t^r),
\]
and hence $\gamma_{\mathbf{u}}(t)\to0$ as $t\to 0^{+}$.

Choose $\tau>0$ with $a:=\gamma_{\mathbf{u}}(\tau)>0$. By the two limits, there exist
$0<\delta<\tau<R$ such that
\[
 \gamma_{\mathbf{u}}(t)<\frac a2
 \quad\text{if }0<t\leq\delta
 \qquad\text{or if }t\geq R.
\]
The continuous function $\gamma_{\mathbf{u}}$ attains a maximum on the compact set
$[\delta,R]$. Since its value at $\tau$ is $a$, this maximum exceeds all
values outside that set and is therefore global. Moreover, it is not attained at the
endpoints $\delta$ or $R$; thus there exists $t_{\mathbf{u}}\in(\delta,R)$ with
$\gamma_{\mathbf{u}}'(t_{\mathbf{u}})=0$, that is, $t_{\mathbf{u}}\mathbf{u}\in\mathcal N$.

We now prove uniqueness. Suppose that $0<s<t$ and $s\mathbf{u},t\mathbf{u}\in\mathcal N$. Dividing the two Nehari identities by $s^q$ and $t^q$, respectively, gives
\[
\int_M a(t^p|\nabla^{\mathbf{E}}\mathbf{u}|_{\mathbf{g}\otimes \mathbf{h}_{\mathbf{E}}}^p)t^{p-q}
|\nabla^{\mathbf{E}}\mathbf{u}|_{\mathbf{g}\otimes \mathbf{h}_{\mathbf{E}}}^p\,d\lambda_{\mathbf{g}}
=
\int_M\frac{\langle f_x(t\mathbf{u}),t\mathbf{u}\rangle_{\mathbf{h}_{\mathbf{E}}}}{t^q}\,d\lambda_{\mathbf{g}}
\]
and the analogous identity with $s$.

If $|\nabla^{\mathbf{E}}\mathbf{u}(x)|_{\mathbf{g}\otimes \mathbf{h}_{\mathbf{E}}}>0$, apply
\textup{(\ref{hip:quasilineal-a2})} to the numbers
$s|\nabla^{\mathbf{E}}\mathbf{u}(x)|_{\mathbf{g}\otimes \mathbf{h}_{\mathbf{E}}}$ and
$t|\nabla^{\mathbf{E}}\mathbf{u}(x)|_{\mathbf{g}\otimes \mathbf{h}_{\mathbf{E}}}$. Since $s<t$, we obtain
\[
a(t^p|\nabla^{\mathbf{E}}\mathbf{u}|_{\mathbf{g}\otimes \mathbf{h}_{\mathbf{E}}}^p)t^{p-q}
\leq
a(s^p|\nabla^{\mathbf{E}}\mathbf{u}|_{\mathbf{g}\otimes \mathbf{h}_{\mathbf{E}}}^p)s^{p-q}.
\]
The inequality also holds when
$|\nabla^{\mathbf{E}}\mathbf{u}(x)|_{\mathbf{g}\otimes \mathbf{h}_{\mathbf{E}}}=0$. Thus the principal term
corresponding to $t$ is less than or equal to the one corresponding to $s$.

By contrast, \textup{(\ref{hip:quasilineal-f4})} implies that, for each $x$ such that $\mathbf{u}(x)\neq0_x$,
\[
\frac{\langle f_x(t\mathbf{u}),t\mathbf{u}\rangle_{\mathbf{h}_{\mathbf{E}}}}{t^q}
>
\frac{\langle f_x(s\mathbf{u}),s\mathbf{u}\rangle_{\mathbf{h}_{\mathbf{E}}}}{s^q}.
\]
The set $\{x\in M\mid \mathbf{u}(x)\neq0_x\}$ has positive measure, since $\mathbf{u}\neq0$ in $L^q(M,\mathbf{E})$. Upon integration, the potential term corresponding to $t$ is strictly greater than the one corresponding to $s$. This contradicts the two Nehari identities. Therefore $t_{\mathbf{u}}$ is unique.

Finally, define $H_{\mathbf{u}}(t):=\displaystyle\frac{DI(t\mathbf{u})(t\mathbf{u})}{t^q}$. The preceding argument, applied to any $0<s<t$, shows that $H_{\mathbf{u}}(t)<H_{\mathbf{u}}(s)$. Since $H_{\mathbf{u}}(t_{\mathbf{u}})=0$, we have $H_{\mathbf{u}}(t)>0$ if $0<t<t_{\mathbf{u}}$ and $H_{\mathbf{u}}(t)<0$ if $t>t_{\mathbf{u}}$. Because $\gamma_{\mathbf{u}}'(t)=t^{q-1}H_{\mathbf{u}}(t)$, the function $\gamma_{\mathbf{u}}$ is strictly increasing on $(0,t_{\mathbf{u}})$ and strictly decreasing on $(t_{\mathbf{u}},+\infty)$.
\end{proof}

\begin{lemma}[Positivity of the energy on $\mathcal N$]
\label{lem:energia-positiva-nehari-quasilineal}
Let $M$ be the compact Riemannian manifold with nonempty boundary fixed in this
section, and let $\mathbf{E}\to M$ be the vector bundle specified there.
For every $\mathbf{u}\in\mathcal N$, we have $I(\mathbf{u})>0$.
\end{lemma}

\begin{proof}
Since $DI(\mathbf{u})(\mathbf{u})=0$, Lemma~\ref{lem:desigualdad-A-a-quasilineal} and \textup{(\ref{hip:quasilineal-f3})} give
\[
\begin{aligned}
I(\mathbf{u})
&=I(\mathbf{u})-\frac{1}{q}DI(\mathbf{u})(\mathbf{u})\\
&=\int_M\left[\frac{1}{p}A(|\nabla^{\mathbf{E}}\mathbf{u}|_{\mathbf{g}\otimes \mathbf{h}_{\mathbf{E}}}^p)
-\frac{1}{q}a(|\nabla^{\mathbf{E}}\mathbf{u}|_{\mathbf{g}\otimes \mathbf{h}_{\mathbf{E}}}^p)
|\nabla^{\mathbf{E}}\mathbf{u}|_{\mathbf{g}\otimes \mathbf{h}_{\mathbf{E}}}^p\right]d\lambda_{\mathbf{g}}\\
&\quad+
\int_M\left[\frac{1}{q}\langle f_x(\mathbf{u}),\mathbf{u}\rangle_{\mathbf{h}_{\mathbf{E}}}
-\mathcal F_x(\mathbf{u})\right]d\lambda_{\mathbf{g}}\\
&\geq
\left(\frac{\theta}{q}-1\right)
\int_M\mathcal F_x(\mathbf{u})\,d\lambda_{\mathbf{g}}>0.
\end{aligned}
\]
The last integral is positive because $\mathbf{u}\neq0$ and \textup{(\ref{hip:quasilineal-f3})} implies $\mathcal F_x(\mathbf{u}(x))>0$ when $\mathbf{u}(x)\neq0_x$.
\end{proof}

\begin{lemma}[Separation from the origin]
\label{lem:nehari-separada-origen-quasilineal}
Let $M$ be the compact Riemannian manifold with nonempty boundary fixed in this
section, and let $\mathbf{E}\to M$ be the vector bundle specified there.
There exists $\rho>0$ such that $\|\mathbf{u}\|_X\geq\rho$ for every $\mathbf{u}\in\mathcal N$.
\end{lemma}

\begin{proof}
Let $\mathbf{u}\in\mathcal N$. By \eqref{eq:identidad-nehari-quasilineal}, \textup{(\ref{hip:quasilineal-a1})}, and \eqref{eq:estimacion-f-quasilineal},
\[
k_1\|\nabla^{\mathbf{E}}\mathbf{u}\|_{L^p(M,T^*M\otimes \mathbf{E})}^p+k_2\|\mathbf{u}\|_X^q
\leq
\varepsilon\|\mathbf{u}\|_{L^p(M,\mathbf{E})}^p+C_\varepsilon\|\mathbf{u}\|_{L^r(M,\mathbf{E})}^r.
\]
Let $C_p,C_r>0$ be constants such that
\[
\|\mathbf{u}\|_{L^p(M,\mathbf{E})}\leq C_p\|\nabla^{\mathbf{E}}\mathbf{u}\|_{L^p(M,T^*M\otimes \mathbf{E})},
\qquad
\|\mathbf{u}\|_{L^r(M,\mathbf{E})}\leq C_r\|\mathbf{u}\|_X.
\]
Take $\varepsilon:=\displaystyle\frac{k_1}{2C_p^p}$. Then
\begin{equation}\label{eq:separacion-nehari-intermedia}
\frac{k_1}{2}\|\nabla^{\mathbf{E}}\mathbf{u}\|_{L^p(M,T^*M\otimes \mathbf{E})}^p+k_2\|\mathbf{u}\|_X^q
\leq
C_\varepsilon C_r^r\|\mathbf{u}\|_X^r.
\end{equation}

If $p<q$, hypothesis \textup{(\ref{hip:quasilineal-a1})} ensures that $k_2>0$. Dropping the first term on the left-hand side of \eqref{eq:separacion-nehari-intermedia} gives $k_2\|\mathbf{u}\|_X^q
\leq
C_\varepsilon C_r^r\|\mathbf{u}\|_X^r$. Since $\mathbf{u}\neq0$ and $r>q$, $\|\mathbf{u}\|_X
\geq
\left(\displaystyle\frac{k_2}{C_\varepsilon C_r^r}\right)^{\frac{1}{r-q}}$.

If $p=q$, we have $\|\nabla^{\mathbf{E}}\mathbf{u}\|_{L^p(M,T^*M\otimes \mathbf{E})}=\|\mathbf{u}\|_X$. Dropping the term containing $k_2$ gives $\displaystyle\frac{k_1}{2}\|\mathbf{u}\|_X^q
\leq
C_\varepsilon C_r^r\|\mathbf{u}\|_X^r$, and hence $\|\mathbf{u}\|_X
\geq
\left(\displaystyle\frac{k_1}{2C_\varepsilon C_r^r}\right)^{\frac{1}{r-q}}$.

Taking $\rho$ to be the smaller of the constants corresponding to the possible cases completes the proof.
\end{proof}

\begin{lemma}[Boundedness of Nehari sequences]
\label{lem:sucesiones-nehari-acotadas-quasilineal}
Let $M$ be the compact Riemannian manifold with nonempty boundary fixed in this
section, and let $\mathbf{E}\to M$ be the vector bundle specified there.
If $(\mathbf{u}_j)_{j\in\mathbb N}\subseteq\mathcal N$ and $(I(\mathbf{u}_j))_{j\in\mathbb N}$ is bounded above, then $(\mathbf{u}_j)$ is bounded in $X$.
\end{lemma}

\begin{proof}
For each $\mathbf{u}\in\mathcal N$, using \textup{(\ref{hip:quasilineal-f3})}, Lemma~\ref{lem:desigualdad-A-a-quasilineal}, and \textup{(\ref{hip:quasilineal-a1})}, we obtain
\[
I(\mathbf{u})=I(\mathbf{u})-\frac{1}{\theta} DI(\mathbf{u})(\mathbf{u})\geq
\int_M\left[
\frac{1}{p}A(|\nabla^{\mathbf{E}}\mathbf{u}|_{\mathbf{g}\otimes \mathbf{h}_{\mathbf{E}}}^p)
-\frac{1}{\theta}a(|\nabla^{\mathbf{E}}\mathbf{u}|_{\mathbf{g}\otimes \mathbf{h}_{\mathbf{E}}}^p)
|\nabla^{\mathbf{E}}\mathbf{u}|_{\mathbf{g}\otimes \mathbf{h}_{\mathbf{E}}}^p
\right]d\lambda_{\mathbf{g}}
\]
\[
\geq
\left(\frac{1}{q}-\frac{1}{\theta}\right)
\int_M a(|\nabla^{\mathbf{E}}\mathbf{u}|_{\mathbf{g}\otimes \mathbf{h}_{\mathbf{E}}}^p)
|\nabla^{\mathbf{E}}\mathbf{u}|_{\mathbf{g}\otimes \mathbf{h}_{\mathbf{E}}}^p\,d\lambda_{\mathbf{g}}\geq
\left(\frac{1}{q}-\frac{1}{\theta}\right)
\left(k_1\|\nabla^{\mathbf{E}}\mathbf{u}\|_{L^p(M,T^*M\otimes \mathbf{E})}^p+k_2\|\mathbf{u}\|_X^q\right).
\]
If $p<q$, then $k_2>0$, and $I(\mathbf{u})
\geq
\left(\displaystyle\frac{1}{q}-\displaystyle\frac{1}{\theta}\right)k_2\|\mathbf{u}\|_X^q$.

If $p=q$, then $\|\nabla^{\mathbf{E}}\mathbf{u}\|_{L^p(M,T^*M\otimes \mathbf{E})}=\|\mathbf{u}\|_X$ and $I(\mathbf{u})
\geq
\left(\displaystyle\frac{1}{q}-\displaystyle\frac{1}{\theta}\right)k_1\|\mathbf{u}\|_X^q$.

In both cases there exists $c_0>0$, independent of $\mathbf{u}\in\mathcal N$, such that $I(\mathbf{u})\geq c_0\|\mathbf{u}\|_X^q$.

The upper bound for $I(\mathbf{u}_j)$ then implies boundedness of $\|\mathbf{u}_j\|_X$.
\end{proof}

\begin{lemma}[Convergence of the nonlinear terms]
\label{lem:convergencia-no-lineal-quasilineal}
Let $M$ be the compact Riemannian manifold with nonempty boundary fixed in this
section, and let $\mathbf{E}\to M$ be the vector bundle specified there.
Let $(\mathbf{u}_j)_{j\in\mathbb N}\subseteq L^r(M,\mathbf{E})$ satisfy $\mathbf{u}_j\to \mathbf{u}$ in $L^r(M,\mathbf{E})$. Then $f(\mathbf{u}_j)\to f(\mathbf{u})
\quad\text{in }L^{r'}(M,\mathbf{E})$, $\displaystyle\int_M\langle f_x(\mathbf{u}_j),\mathbf{u}_j\rangle_{\mathbf{h}_{\mathbf{E}}}\,d\lambda_{\mathbf{g}}
\longrightarrow
\int_M\langle f_x(\mathbf{u}),\mathbf{u}\rangle_{\mathbf{h}_{\mathbf{E}}}\,d\lambda_{\mathbf{g}}$, and $\displaystyle\int_M\mathcal F_x(\mathbf{u}_j)\,d\lambda_{\mathbf{g}}
\longrightarrow
\int_M\mathcal F_x(\mathbf{u})\,d\lambda_{\mathbf{g}}$.

\end{lemma}

\begin{proof}
The first convergence follows from Corollary~\ref{cor:superposicion-exponentes-conjugados-haces}, applied to the bundle map $f\colon \mathbf{E}\longrightarrow \mathbf{E}$ with $m=r$.

Since $\mathbf{u}_j\to \mathbf{u}$ in $L^r(M,\mathbf{E})$, the sequence $(\mathbf{u}_j)$ is bounded in $L^r(M,\mathbf{E})$. By Hölder's inequality (Proposition~\ref{desigualdad de holder}),
\[
\begin{aligned}
&\left|
\int_M\langle f_x(\mathbf{u}_j),\mathbf{u}_j\rangle_{\mathbf{h}_{\mathbf{E}}}\,d\lambda_{\mathbf{g}}
-
\int_M\langle f_x(\mathbf{u}),\mathbf{u}\rangle_{\mathbf{h}_{\mathbf{E}}}\,d\lambda_{\mathbf{g}}
\right|\\
&\qquad\leq
\|f(\mathbf{u}_j)-f(\mathbf{u})\|_{L^{r'}(M,\mathbf{E})}\|\mathbf{u}_j\|_{L^r(M,\mathbf{E})}
+
\|f(\mathbf{u})\|_{L^{r'}(M,\mathbf{E})}\|\mathbf{u}_j-\mathbf{u}\|_{L^r(M,\mathbf{E})}.
\end{aligned}
\]
and the right-hand side converges to zero.

For the potential, Theorem~\ref{teo:b5-fundamental-calculo-riemann-banach}, applied in each fiber, gives
\[
\mathcal F_x(\mathbf{u}_j)-\mathcal F_x(\mathbf{u})
=
\int_0^1\langle f_x(\mathbf{u}+s(\mathbf{u}_j-\mathbf{u})),\mathbf{u}_j-\mathbf{u}\rangle_{\mathbf{h}_{\mathbf{E}}}\,ds.
\]
By Fubini's theorem for Bochner integrals~\ref{teo:b5-fubini-bochner}, applied in the scalar case, and Hölder's inequality (Proposition~\ref{desigualdad de holder}),
\[
\begin{split}
\left|
\int_M\mathcal F_x(\mathbf{u}_j)\,d\lambda_{\mathbf{g}}
-
\int_M\mathcal F_x(\mathbf{u})\,d\lambda_{\mathbf{g}}
\right|
&\leq
\|\mathbf{u}_j-\mathbf{u}\|_{L^r(M,\mathbf{E})}
\int_0^1\|f(\mathbf{u}+s(\mathbf{u}_j-\mathbf{u}))\|_{L^{r'}(M,\mathbf{E})}\,ds.
\end{split}
\]
The family $\mathbf{u}+s(\mathbf{u}_j-\mathbf{u})$ is bounded in $L^r(M,\mathbf{E})$ uniformly with respect to $j$ and $s\in[0,1]$. Remark~\ref{rem: nemytski-variedad-compacta-haz} shows that the norms of $f(\mathbf{u}+s(\mathbf{u}_j-\mathbf{u}))$ in $L^{r'}(M,\mathbf{E})$ are uniformly bounded. Since $\|\mathbf{u}_j-\mathbf{u}\|_{L^r(M,\mathbf{E})}\to0$, the difference of the potentials converges to zero.
\end{proof}

\subsection{Radial reduction and existence of a minimizer}

Proposition~\ref{prop:proyeccion-nehari-quasilineal} allows us to define $m\colon S_X\longrightarrow\mathcal N$ by $m(\mathbf{w})=t_{\mathbf{w}}\mathbf{w}$. To apply Theorem~\ref{teo:reduccion-radial-szulkin-weth}, it remains to verify the uniform bounds required there.

\begin{proposition}[Hypotheses for radial reduction]
\label{prop:hipotesis-szulkin-weth-quasilineal}
Let $M$ be the compact Riemannian manifold with nonempty boundary fixed in this
section, and let $\mathbf{E}\to M$ be the vector bundle specified there.
The space $X$ and the functional $I$ satisfy hypotheses \textup{(\ref{hip:szulkin-A1})}--\textup{(\ref{hip:szulkin-A3})} for the radial reduction of Szulkin and Weth.
\end{proposition}

\begin{proof}
Consider the linear map
\[
 \mathcal G\colon X\longrightarrow L^q(M,T^*M\otimes\mathbf E),
 \qquad
 \mathcal G(\mathbf u)=\nabla^{\mathbf E}\mathbf u.
\]
The Poincaré inequality in
Theorem~\ref{teo:poincare-W0-haces} shows that $\mathcal G$ is
injective. Moreover, by the definition of the norm fixed in this section,
\[
 \|\mathcal G(\mathbf u)\|_{L^q(M,T^*M\otimes\mathbf E)}
 =\|\mathbf u\|_X.
\]
The norm $\|\cdot\|_X$ is equivalent to the norm on
$W_0^{1,q}(M,\mathbf E)$; since the latter space is closed in
$W^{1,q}(M,\mathbf E)$ and complete, $(X,\|\cdot\|_X)$ is a
Banach space. Consequently, the image $\mathcal G(X)$ is closed in
$L^q(M,T^*M\otimes\mathbf E)$. The proof of
Proposition~\ref{prop: Lp(E) es uniformemente convexo} is pointwise in the
fibers and uses only that they are Hilbert spaces; it therefore applies
without change on the manifold with boundary and shows that
$L^q(M,T^*M\otimes\mathbf E)$ is uniformly convex. This property is
inherited by closed subspaces and preserved under isometries. Thus
$X$ is uniformly convex.

To verify \textup{(\ref{hip:szulkin-A1})}, take the normalization function $\varphi(t):=t^{q-1}$. The associated functional is
\[
\psi_X(\mathbf{u}):=\int_0^{\|\mathbf{u}\|_X}t^{q-1}\,dt
=\frac{1}{q}\|\mathbf{u}\|_X^q
=\frac{1}{q}\int_M|\nabla^{\mathbf{E}}\mathbf{u}|_{\mathbf{g}\otimes \mathbf{h}_{\mathbf{E}}}^q\,d\lambda_{\mathbf{g}}.
\]
The functional $\psi_X$ belongs to $C^1(X,\mathbb R)$. This is the case $p=q$, $a\equiv1$, and $A(s)=s$ of the calculation for $J_1$ in the proof of Proposition~\ref{prop:energia-quasilineal-C1}. In this specialization,
\[
D\psi_X(\mathbf{u})(\mathbf{v})
=
\int_M
|\nabla^{\mathbf{E}}\mathbf{u}|_{\mathbf{g}\otimes \mathbf{h}_{\mathbf{E}}}^{q-2}
\langle\nabla^{\mathbf{E}}\mathbf{u},\nabla^{\mathbf{E}}\mathbf{v}\rangle_{\mathbf{g}\otimes \mathbf{h}_{\mathbf{E}}}\,d\lambda_{\mathbf{g}}.
\]
Hölder's inequality (Proposition~\ref{desigualdad de holder}) implies $|D\psi_X(\mathbf{u})(\mathbf{v})|
\leq
\|\mathbf{u}\|_X^{q-1}\|\mathbf{v}\|_X$, so the map $\mathcal J_X:=D\psi_X\colon X\longrightarrow X'$ sends bounded sets to bounded sets. If $\mathbf{w}\in S_X$, then
\[
\mathcal J_X(\mathbf{w})(\mathbf{w})
=
\int_M|\nabla^{\mathbf{E}}\mathbf{w}|_{\mathbf{g}\otimes \mathbf{h}_{\mathbf{E}}}^q\,d\lambda_{\mathbf{g}}
=
\|\mathbf{w}\|_X^q
=1.
\]
Thus \textup{(\ref{hip:szulkin-A1})} holds.

Proposition~\ref{prop:proyeccion-nehari-quasilineal} proves that, for each $\mathbf{u}\in X\setminus\{0\}$, there exists a unique $t_{\mathbf{u}}>0$ such that $t_{\mathbf{u}}\mathbf{u}\in\mathcal N$, and that the function $t\mapsto I(t\mathbf{u})$ is strictly increasing on $(0,t_{\mathbf{u}})$ and strictly decreasing on $(t_{\mathbf{u}},+\infty)$. This is precisely hypothesis \textup{(\ref{hip:szulkin-A2})}.

We now verify \textup{(\ref{hip:szulkin-A3})}. If $\mathbf{w}\in S_X$, then $t_{\mathbf{w}}\mathbf{w}\in\mathcal N$. By Lemma~\ref{lem:nehari-separada-origen-quasilineal}, there exists $\rho>0$, independent of $\mathbf{w}$, such that $t_{\mathbf{w}}=\|t_{\mathbf{w}}\mathbf{w}\|_X\geq\rho$. This provides the uniform lower bound.

Now let $K\subseteq S_X$ be nonempty and compact; for $K=\varnothing$ the upper bound is vacuous. The function $\mathbf{w}\longmapsto\|\mathbf{w}\|_{L^\theta(M,\mathbf{E})}^\theta$ is continuous on $X$ and positive on $S_X$, since $\|\mathbf{w}\|_{L^\theta(M,\mathbf{E})}=0$ would imply $\mathbf{w}=0$. By compactness, $\displaystyle \mu_K:=\min_{\mathbf{w}\in K}\|\mathbf{w}\|_{L^\theta(M,\mathbf{E})}^\theta>0$.

By the upper bound for $A$ and Lemma~\ref{lem:cota-inferior-potencial-quasilineal}, there exists $C_K>0$ such that, for every $\mathbf{w}\in K$ and every $t>0$,
\[
I(t\mathbf{w})
\leq
C_K(t^p+t^q)-D\mu_Kt^\theta+D'\lambda_{\mathbf{g}}(M).
\]
Since $\theta>q\geq p$, the right-hand side converges to $-\infty$ as $t\to+\infty$, uniformly with respect to $\mathbf{w}\in K$. There exists $R_K>0$ such that $I(t\mathbf{w})<0$ for every $\mathbf{w}\in K$ and every $t\geq R_K$. However, $I(t_{\mathbf{w}}\mathbf{w})>0$ by Lemma~\ref{lem:energia-positiva-nehari-quasilineal}. Thus $t_{\mathbf{w}}<R_K$ for every $\mathbf{w}\in K$, proving \textup{(\ref{hip:szulkin-A3})}.
\end{proof}

\begin{proposition}[Existence of a minimizer on $\mathcal N$]
\label{prop:minimizador-nehari-quasilineal}
Let $M$ be the compact Riemannian manifold with nonempty boundary fixed in this
section, and let $\mathbf{E}\to M$ be the vector bundle specified there.
There exists $\mathbf{u}_0\in\mathcal N$ such that
\[
I(\mathbf{u}_0)=\inf_{\mathbf{u}\in\mathcal N}I(\mathbf{u}).
\]
\end{proposition}

\begin{proof}
By Proposition~\ref{prop:proyeccion-nehari-quasilineal}, every nonzero
section determines a point of $\mathcal N$, so $\mathcal N\neq\varnothing$.
Set $c:=\displaystyle\inf_{\mathbf u\in\mathcal N}I(\mathbf u)$.
Lemma~\ref{lem:sucesiones-nehari-acotadas-quasilineal} gives $c_0>0$
with $I(\mathbf u)\geq c_0\|\mathbf u\|_X^q$ on $\mathcal N$.
If $\rho>0$ is the constant in
Lemma~\ref{lem:nehari-separada-origen-quasilineal}, then
\[
 0<c_0\rho^q\leq c\leq I(\mathbf v)<\infty
 \qquad(\mathbf v\in\mathcal N).
\]
We can therefore choose $(\mathbf u_j)_{j\in\mathbb N}\subseteq\mathcal N$
with $c\leq I(\mathbf u_j)<c+1/j$.
Lemma~\ref{lem:sucesiones-nehari-acotadas-quasilineal} bounds this sequence
in $X$. Reflexivity gives a subsequence and $\mathbf u\in X$ with
$\mathbf u_j\rightharpoonup\mathbf u$ in $X$.

Since $q<r<q^*$, Corollary~\ref{cor:encajes-W0-haces-frontera} yields a
further subsequence converging in $L^r(M,\mathbf E)$.
Its limit is $\mathbf u$: the continuous embedding carries weak
convergence in $X$ to weak convergence in $L^r$, and a strong limit
has that same weak limit. Relabeling this subsequence, we have
\[
 \mathbf u_j\rightharpoonup\mathbf u\text{ in }X,
 \qquad \mathbf u_j\to\mathbf u\text{ in }L^r(M,\mathbf E).
\]
Moreover, $p<r$ and the volume is finite, so Hölder gives
\[
 \|\mathbf u_j-\mathbf u\|_{L^p(M,\mathbf E)}
 \leq\lambda_{\mathbf g}(M)^{1/p-1/r}
       \|\mathbf u_j-\mathbf u\|_{L^r(M,\mathbf E)}\longrightarrow0.
\]

We shall prove that $\mathbf{u}\neq0$. Suppose otherwise. By Lemma~\ref{lem:convergencia-no-lineal-quasilineal}, $\displaystyle\int_M\langle f_x(\mathbf{u}_j),\mathbf{u}_j\rangle_{\mathbf{h}_{\mathbf{E}}}\,d\lambda_{\mathbf{g}}\longrightarrow0$.

Fix $C_p>0$ such that $\|\mathbf{v}\|_{L^p(M,\mathbf{E})}^p\leq C_p\|\nabla^{\mathbf{E}}\mathbf{v}\|_{L^p(M,T^*M\otimes \mathbf{E})}^p$ for every $\mathbf{v}\in X$, and choose $\varepsilon:=\displaystyle\frac{k_1}{2C_p}$. By the Nehari identity, \textup{(\ref{hip:quasilineal-a1})}, and \eqref{eq:estimacion-f-quasilineal},
\[
\begin{aligned}
&k_1\|\nabla^{\mathbf{E}}\mathbf{u}_j\|_{L^p(M,T^*M\otimes \mathbf{E})}^p+k_2\|\mathbf{u}_j\|_X^q\\
&\qquad\leq
\int_M\langle f_x(\mathbf{u}_j),\mathbf{u}_j\rangle_{\mathbf{h}_{\mathbf{E}}}\,d\lambda_{\mathbf{g}}\\
&\qquad\leq
\varepsilon\|\mathbf{u}_j\|_{L^p(M,\mathbf{E})}^p+C_\varepsilon\|\mathbf{u}_j\|_{L^r(M,\mathbf{E})}^r\\
&\qquad\leq
\frac{k_1}{2}\|\nabla^{\mathbf{E}}\mathbf{u}_j\|_{L^p(M,T^*M\otimes \mathbf{E})}^p
+C_\varepsilon\|\mathbf{u}_j\|_{L^r(M,\mathbf{E})}^r.
\end{aligned}
\]
Thus
\[
\frac{k_1}{2}\|\nabla^{\mathbf{E}}\mathbf{u}_j\|_{L^p(M,T^*M\otimes \mathbf{E})}^p+k_2\|\mathbf{u}_j\|_X^q
\leq
C_\varepsilon\|\mathbf{u}_j\|_{L^r(M,\mathbf{E})}^r.
\]
If $p<q$, we have $k_2>0$ and the right-hand side converges to zero, so $\|\mathbf{u}_j\|_X\to0$. If $p=q$, the first term on the left-hand side is $\displaystyle\frac{k_1}{2}\|\mathbf{u}_j\|_X^q$, giving the same conclusion. This contradicts $\|\mathbf{u}_j\|_X\geq\rho$. Hence $\mathbf{u}\neq0$.

By Proposition~\ref{prop:proyeccion-nehari-quasilineal}, there exists a unique $t_{\mathbf{u}}>0$ such that $t_{\mathbf{u}}\mathbf{u}\in\mathcal N$. Since $\mathbf{u}_j\in\mathcal N$ is the maximum along its own positive ray, $I(t_{\mathbf{u}}\mathbf{u}_j)\leq I(\mathbf{u}_j)$ for every $j$.

Set $\phi(s):=A(s^p)$ for $s\geq0$. By
\textup{(\ref{hip:quasilineal-a3})}, $\phi$ is convex; since $a\geq0$,
$A$ and $\phi$ are also nondecreasing. The triangle inequality and these
two properties give, for $0\leq\lambda\leq1$,
\[
\begin{aligned}
A\bigl(|\lambda\eta+(1-\lambda)\zeta|_{\mathbf{g}\otimes \mathbf{h}_{\mathbf{E}}}^p\bigr)
&=\phi\bigl(|\lambda\eta+(1-\lambda)\zeta|_{\mathbf{g}\otimes \mathbf{h}_{\mathbf{E}}}\bigr)\\
&\leq\phi\bigl(
\lambda|\eta|_{\mathbf{g}\otimes \mathbf{h}_{\mathbf{E}}}
+(1-\lambda)|\zeta|_{\mathbf{g}\otimes \mathbf{h}_{\mathbf{E}}}\bigr)\\
&\leq\lambda A(|\eta|_{\mathbf{g}\otimes \mathbf{h}_{\mathbf{E}}}^p)
+(1-\lambda)A(|\zeta|_{\mathbf{g}\otimes \mathbf{h}_{\mathbf{E}}}^p).
\end{aligned}
\]
Thus $\eta\mapsto p^{-1}A(|\eta|_{\mathbf{g}\otimes \mathbf{h}_{\mathbf{E}}}^p)$ is convex and
continuous. Consequently, for the fixed number $t_{\mathbf{u}}>0$, the functional
\[
G(\mathbf{v}):=\frac{1}{p}\int_M A(t_{\mathbf{u}}^p|\nabla^{\mathbf{E}}\mathbf{v}|_{\mathbf{g}\otimes \mathbf{h}_{\mathbf{E}}}^p)\,d\lambda_{\mathbf{g}}
\]
is convex and continuous on $X$. By Corollary~\ref{cor: funcional convexo semicontinuo es debilmente semicontinuo}, $G$ is weakly lower semicontinuous; therefore,
\[
\frac{1}{p}\int_M
A(t_{\mathbf{u}}^p|\nabla^{\mathbf{E}}\mathbf{u}|_{\mathbf{g}\otimes \mathbf{h}_{\mathbf{E}}}^p)\,d\lambda_{\mathbf{g}}
\leq
\liminf_{j\to\infty}
\frac{1}{p}\int_M
A(t_{\mathbf{u}}^p|\nabla^{\mathbf{E}}\mathbf{u}_j|_{\mathbf{g}\otimes \mathbf{h}_{\mathbf{E}}}^p)\,d\lambda_{\mathbf{g}}.
\]
Moreover, $t_{\mathbf{u}}\mathbf{u}_j\to t_{\mathbf{u}}\mathbf{u}$ in $L^r(M,\mathbf{E})$, so Lemma~\ref{lem:convergencia-no-lineal-quasilineal} implies
\[
\int_M\mathcal F_x(t_{\mathbf{u}}\mathbf{u}_j)\,d\lambda_{\mathbf{g}}
\longrightarrow
\int_M\mathcal F_x(t_{\mathbf{u}}\mathbf{u})\,d\lambda_{\mathbf{g}}.
\]
Consequently, $\displaystyle I(t_{\mathbf{u}}\mathbf{u})
\leq
\liminf_{j\to\infty}I(t_{\mathbf{u}}\mathbf{u}_j)
\leq
\lim_{j\to\infty}I(\mathbf{u}_j)
=c$.

Since $t_{\mathbf{u}}\mathbf{u}\in\mathcal N$, the definition of $c$ gives $c\leq I(t_{\mathbf{u}}\mathbf{u})$. Thus $I(t_{\mathbf{u}}\mathbf{u})=c$. Taking $\mathbf{u}_0:=t_{\mathbf{u}}\mathbf{u}$ completes the proof.
\end{proof}

\begin{theorem}[Existence of a ground state]
\label{teo:estado-fundamental-quasilineal-haz}
Let $M$ be the compact Riemannian manifold with nonempty boundary fixed in this
section, and let $\mathbf{E}\to M$ be the vector bundle specified there.
Under hypotheses \textup{(\ref{hip:quasilineal-a1})}--\textup{(\ref{hip:quasilineal-a3})} and \textup{(\ref{hip:quasilineal-f1})}--\textup{(\ref{hip:quasilineal-f4})}, problem \eqref{eq:problema-quasilineal-haz} has a ground state solution: there exists $\mathbf{u}_0\in X\setminus\{0\}$ such that $DI(\mathbf{u}_0)=0$ and
\[
I(\mathbf{u}_0)
=
\inf\{I(\mathbf{u})\mid \mathbf{u}\in X\setminus\{0\},\ DI(\mathbf{u})=0\}.
\]
\end{theorem}

\begin{proof}
By Proposition~\ref{prop:minimizador-nehari-quasilineal}, there exists
$\mathbf{u}_0\in\mathcal N$ minimizing $I$ on $\mathcal N$.
Propositions~\ref{prop:proyeccion-nehari-quasilineal} and
\ref{prop:hipotesis-szulkin-weth-quasilineal} allow us to apply
Theorem~\ref{teo:reduccion-radial-szulkin-weth}. Let
$\mathbf{w}_0:=m^{-1}(\mathbf{u}_0)\in S_X$. Since $m$ is a homeomorphism and
$\Psi=I\circ m$, the point $\mathbf{w}_0$ minimizes $\Psi$ on the submanifold
$S_X$. Fermat's theorem on manifolds gives $D\Psi(\mathbf{w}_0)=0$; by
part~(e) of Theorem~\ref{teo:reduccion-radial-szulkin-weth},
$\mathbf{u}_0=m(\mathbf{w}_0)$ is a nontrivial critical point of $I$. By
Proposition~\ref{prop:energia-quasilineal-C1}, $\mathbf{u}_0$ is a weak solution
of \eqref{eq:problema-quasilineal-haz}.

Every nontrivial critical point $\mathbf{u}$ of $I$ satisfies $DI(\mathbf{u})(\mathbf{u})=0$ and therefore belongs to $\mathcal N$. It follows that
\[
I(\mathbf{u}_0)
=
\inf_{\mathbf{u}\in\mathcal N}I(\mathbf{u})
=
\inf\{I(\mathbf{u})\mid \mathbf{u}\neq0,\ DI(\mathbf{u})=0\}.
\]
\end{proof}

\subsection{Examples of operators covered by the result}

The following examples show that the hypotheses on $a$ encompass several familiar geometric operators. In all of them, we may take
\[
\mathcal F_x(\xi):=\frac{1}{r}|\xi|_{\mathbf{h}_{\mathbf{E}}}^r,
\qquad
f_x(\xi)=|\xi|_{\mathbf{h}_{\mathbf{E}}}^{r-2}\xi,
\qquad
r\in(q,q^*).
\]
These conditions are verified directly, since $\displaystyle\frac{|f_x(\xi)|_{\mathbf{h}_{\mathbf{E}}}}{|\xi|_{\mathbf{h}_{\mathbf{E}}}^{p-1}}
=
|\xi|_{\mathbf{h}_{\mathbf{E}}}^{r-p}
\longrightarrow0$ uniformly with respect to $x$, because $r>q\geq p$; moreover, $|f_x(\xi)|_{\mathbf{h}_{\mathbf{E}}}=|\xi|_{\mathbf{h}_{\mathbf{E}}}^{r-1}$, $\langle f_x(\xi),\xi\rangle_{\mathbf{h}_{\mathbf{E}}}=|\xi|_{\mathbf{h}_{\mathbf{E}}}^r=r\mathcal F_x(\xi)$, and, for $\xi\neq0_x$, $\displaystyle\frac{\langle f_x(t\xi),t\xi\rangle_{\mathbf{h}_{\mathbf{E}}}}{t^q}
=
t^{r-q}|\xi|_{\mathbf{h}_{\mathbf{E}}}^r$, which is strictly increasing. Thus \textup{(\ref{hip:quasilineal-f1})}--\textup{(\ref{hip:quasilineal-f4})} hold for any $\theta\in(q,r]$.

\begin{example}[Bochner Laplacian of type $p$]
Suppose that $p=q$ and $a(s):=1$. Then $a(t^p)t^p=t^p$, so \textup{(\ref{hip:quasilineal-a1})} holds with $k_1=k_3=1$ and $k_2=k_4=0$. Since $\displaystyle\frac{a(t^p)}{t^{q-p}}=1$, \textup{(\ref{hip:quasilineal-a2})} holds. Finally, $A(s)=s$ and $A(t^p)=t^p$ is convex, so \textup{(\ref{hip:quasilineal-a3})} also holds. The principal term in \eqref{eq:problema-quasilineal-haz} reduces to $(\nabla^{\mathbf{E}})^*\bigl(|\nabla^{\mathbf{E}}\mathbf{u}|_{\mathbf{g}\otimes \mathbf{h}_{\mathbf{E}}}^{p-2}\nabla^{\mathbf{E}}\mathbf{u}\bigr)$.

\end{example}

\begin{example}[Bochner operator with $(p,q)$ growth]
Let $a(s):=1+s^{\frac{q-p}{p}}$. Then $a(t^p)t^p=t^p+t^q$, so \textup{(\ref{hip:quasilineal-a1})} holds with $k_1=k_2=k_3=k_4=1$ when $p<q$. If $p=q$, then $a(s)=2$ for every $s\geq0$, and \textup{(\ref{hip:quasilineal-a1})} holds, for example, with $k_1=k_3=2$ and $k_2=k_4=0$. For $p<q$, $\displaystyle\frac{a(t^p)}{t^{q-p}}
=
1+t^{p-q}$, which is strictly decreasing; for $p=q$ it is constant. Moreover,
\[
A(s)=s+\frac{p}{q}s^{\frac{q}{p}},
\qquad
A(t^p)=t^p+\frac{p}{q}t^q,
\]
and the latter function is convex. Thus
\textup{(\ref{hip:quasilineal-a1})}--\textup{(\ref{hip:quasilineal-a3})}
hold, and the principal term is the sum
$(\nabla^{\mathbf{E}})^*\bigl(
|\nabla^{\mathbf{E}}\mathbf{u}|_{\mathbf{g}\otimes \mathbf{h}_{\mathbf{E}}}^{p-2}\nabla^{\mathbf{E}}\mathbf{u}\bigr)
+
(\nabla^{\mathbf{E}})^*\bigl(
|\nabla^{\mathbf{E}}\mathbf{u}|_{\mathbf{g}\otimes \mathbf{h}_{\mathbf{E}}}^{q-2}\nabla^{\mathbf{E}}\mathbf{u}\bigr)$.
\end{example}

\begin{example}[Perturbation of mean curvature type]
Suppose that $q=p$ and
\[
a(s):=1+\frac{1}{(1+s)^{\frac{p-1}{p}}}.
\]
Since $1\leq a(s)\leq2$, we have $t^p\leq a(t^p)t^p\leq2t^p$, so \textup{(\ref{hip:quasilineal-a1})} holds with $k_1=1$, $k_2=0$, $k_3=2$, and $k_4=0$. Because $q-p=0$ and $t\mapsto a(t^p)$ is nonincreasing, \textup{(\ref{hip:quasilineal-a2})} holds. Direct integration gives $A(s)=s+p\bigl((1+s)^{\frac{1}{p}}-1\bigr)$, and hence $A(t^p)=t^p+p\bigl((1+t^p)^{\frac{1}{p}}-1\bigr)$.

The function $t\mapsto t^p$ is convex, and $t\mapsto(1+t^p)^{\frac{1}{p}}$ is the restriction to $[0,+\infty)$ of the $\ell^p(\{1,2\})$ norm of the vector $(1,t)$; consequently, $A(t^p)$ is convex and \textup{(\ref{hip:quasilineal-a3})} holds. In this case, the additional term
\[
(\nabla^{\mathbf{E}})^*\left(
\frac{|\nabla^{\mathbf{E}}\mathbf{u}|_{\mathbf{g}\otimes \mathbf{h}_{\mathbf{E}}}^{p-2}\nabla^{\mathbf{E}}\mathbf{u}}
{(1+|\nabla^{\mathbf{E}}\mathbf{u}|_{\mathbf{g}\otimes \mathbf{h}_{\mathbf{E}}}^{p})^{\frac{p-1}{p}}}
\right).
\]
appears.
\end{example}

\begin{remark}
The condition $k_2>0$ is required only when $p<q$. This distinction is necessary in the coercive estimates on $\mathcal N$: if $p<q$, the term $k_2\|\mathbf{u}\|_X^q$ directly controls the norm of $X$; if $p=q$, this control is provided by the term with coefficient $k_1$, because $\|\nabla^{\mathbf{E}}\mathbf{u}\|_{L^p(M,T^*M\otimes \mathbf{E})}=\|\mathbf{u}\|_X$.
\end{remark}

\section{Exercises}

\begin{exercise}
Let $(M,\mathbf{g})$ be a nonempty compact Riemannian manifold with smooth boundary, such that each connected component has nonempty boundary. For the functional in the model case, prove that for each $u\neq0$ there exists a unique $t_u>0$ such that $t_u u\in\mathcal N$, and compute $t_u$ explicitly. Also show that $u\mapsto t_u u$ is a homeomorphism between the unit sphere of $W_0^{1,2}(M)$ and $\mathcal N$.
\end{exercise}

\begin{exercise}
Let $(M,\mathbf{g})$ be a nonempty compact Riemannian manifold with smooth boundary, such that each connected component has nonempty boundary, and let $\mathbf{E}\to M$ be a smooth real vector bundle of positive finite rank, equipped with a bundle metric $\mathbf{h}_{\mathbf{E}}$ and a connection $\nabla^{\mathbf{E}}$ compatible with $\mathbf{h}_{\mathbf{E}}$. Suppose that $a$ satisfies \textup{(\ref{hip:quasilineal-a1})}--\textup{(\ref{hip:quasilineal-a3})} and that $\mathcal F$ satisfies \textup{(\ref{hip:quasilineal-f1})}--\textup{(\ref{hip:quasilineal-f4})}. Prove directly, without using the existence of the minimizer, that for each $\mathbf{u}\in X\setminus\{0\}$ there exist numbers $0<\alpha_{\mathbf{u}}<\beta_{\mathbf{u}}$ such that $I(t\mathbf{u})>0$ if $0<t<\alpha_{\mathbf{u}}$ and $I(t\mathbf{u})<0$ if $t>\beta_{\mathbf{u}}$. As a hint, first obtain explicit constants $C_1(\mathbf{u}),C_2(\mathbf{u})>0$ for which
\[
I(t\mathbf{u})\geq C_1(\mathbf{u})t^p-C_2(\mathbf{u})t^r
\]
when $t>0$, and then use the lower bound of order $\theta$ for the potential. Explain why this argument alone does not establish uniqueness of the critical point of the radial function, and identify the precise role of \textup{(\ref{hip:quasilineal-a2})} and \textup{(\ref{hip:quasilineal-f4})} in the proof of uniqueness.
\end{exercise}

\appendix
\part*{Appendices}
\chapter{Topology}\label{ap:topologia}

This appendix collects the topological notions used in constructing
manifolds, bundles, and quotient spaces, together with the results from
algebraic topology needed for index theory. For the latter results, we
specify the conventions used in the main chapters and give precise
references for the topological proofs themselves.

Before beginning, we state the maximality principle that will be used to
construct separated nets in metric spaces.

\begin{theorem}[Zorn's lemma]\label{teo:lema-zorn}\index{Zorn's lemma}
Let $(P,\preccurlyeq)$ be a nonempty partially ordered set. If every
chain $\mathcal C\subseteq P$ has an upper bound in $P$, then
$P$ contains a maximal element.
\end{theorem}

\begin{remark}
Zorn's lemma is equivalent to the axiom of choice. Besides constructing
maximal families of separated balls, the text uses choice in extending
linearly independent sets to Hamel bases, in well-ordering index sets for
transfinite constructions, and through the Hahn--Banach theorem.
These uses occur, in particular, in the chapters on Banach bundles and
partitions of unity and in Appendix~\ref{ap:analisis-matematico}.
\end{remark}

\section{Topological spaces}
\begin{definition}[Topological space]\label{def:a-topologia-espacio-topologico}\index{topological space@topological space}
A \textit{topological space} is a pair $(X,\tau)$, where $X$ is a set and $\tau$ is a family of subsets of $X$ such that:
\begin{enumerate}[label=(\alph*)]
 \item $\varnothing, X \in \tau$,
 \item every arbitrary union of sets in $\tau$ belongs to $\tau$,
 \item every finite intersection of sets in $\tau$ belongs to $\tau$.
\end{enumerate}
The elements of $\tau$ are called \textit{open subsets} of $X$.
\end{definition}

\begin{definition}[Subspace topology]\label{def:a-topologia-topologia-de-subespacio}\index{subspace topology@subspace topology}
Let $(X,\tau)$ be a topological space and let $A\subseteq X$. The \textit{subspace topology} on $A$ is
\[
\tau_{A}=\{U\cap A\mid U\in \tau\}.
\]
Then $(A,\tau_{A})$ is a topological space called a \textit{subspace} of $X$.
\end{definition}

\begin{definition}[Quotient map and quotient topology]\label{def:a-topologia-mapeo-cociente-y-topologia-cociente}\index{quotient map}\index{quotient topology@quotient topology}
Let $(X,\tau_X)$ be a topological space, let $Y$ be a set, and let
$q\colon X\longrightarrow Y$ be a surjective map. The \textit{quotient topology} on $Y$ induced by $q$ is defined by
\[
\tau_q=\{V\subseteq Y \mid q^{-1}(V)\in\tau_X\}.
\]
With this topology, $q\colon (X,\tau_X)\longrightarrow(Y,\tau_q)$ is
continuous, and $q$ is called a \textit{quotient map}. Equivalently,
a continuous surjective map between topological spaces is a quotient map if
$U\subseteq Y$ is open whenever $q^{-1}(U)$ is open in $X$.
\end{definition}
\begin{definition}[Quotient space]\label{def:a-topologia-espacio-cociente-2}\index{quotient space}
Let $X$ be a topological space and let $\sim$ be an equivalence relation on $X$.
The \textit{quotient space} $X/\sim$ is defined as the set of equivalence classes
\[
X/\sim =\{[x]\mid x\in X\},
\]
equipped with the quotient topology induced by the natural projection
\[
q\colon X\longrightarrow X/\sim,\qquad q(x)=[x].
\]
Thus, $U\subseteq X/\sim$ is open if and only if $q^{-1}(U)$ is open in $X$.
\end{definition}

\begin{example}[Circle]\label{ej:a-topologia-circunferencia}
The circle $\mathbb S^1$ can be described as the quotient
\[
\mathbb S^1\cong[0,1]/\sim,
\]
where $0\sim 1$ and all other points of $[0,1]$ remain distinct.
\end{example}

\begin{example}[Cylinder]\label{ej:a-topologia-cilindro}
The cylinder is obtained as
\[
C\cong\bigl([0,1]\times[0,1]\bigr)/\sim,
\]
where we identify $(0,y)\sim(1,y)$ for every $y\in [0,1]$. In this case the vertical sides are glued together without a twist.
\end{example}

\begin{example}[Möbius strip]\label{ej:a-topologia-banda-de-mobius}
The Möbius strip is obtained similarly as
\[
M\cong\bigl([0,1]\times[0,1]\bigr)/\sim,
\]
where the equivalence relation is
\[
(0,y)\sim(1,1-y),\qquad y\in [0,1].
\]
This corresponds to gluing the vertical sides of the square after giving one of them a half-turn.
\end{example}

\begin{lemma}[Characterization of quotient maps]\label{lem:a-topologia-caracterizacion-cociente}\index{characterization of quotient maps@characterization of quotient maps}
Let $q\colon X\longrightarrow Y$ be a continuous surjective map between topological spaces.
If $q$ is closed, then $q$ is a quotient map; that is, a subset $U\subseteq Y$ is open if and only if $q^{-1}(U)$ is open in $X$.
\end{lemma}

\begin{proof}
Continuity proves one implication. Conversely, if
$q^{-1}(U)$ is open, then $q^{-1}(Y\setminus U)$ is closed. Since
$q$ is closed and surjective,
\[
 Y\setminus U=q\bigl(q^{-1}(Y\setminus U)\bigr)
\]
is closed. Therefore, $U$ is open.
\end{proof}

\section{Covering spaces}
\begin{definition}[Covering map]\label{def:a-topologia-mapeo-cubriente}\index{covering map}
Let $p\colon Y\longrightarrow X$ be a continuous surjective map between topological spaces. We say that $p$ is a \textit{covering map} if for every $x\in X$ there is an open subset $U\subseteq X$ with $x\in U$ such that
\[
p^{-1}(U)=\coprod_{\alpha\in A} V_{\alpha},
\]
where each $V_{\alpha}$ is open in $Y$ and the restriction $p\restriction_{V_{\alpha}}\colon V_{\alpha}\longrightarrow U$ is a homeomorphism.
\end{definition}

\begin{definition}[Universal covering map]\label{def:a-topologia-mapeo-cubriente-universal}\index{universal covering map}
A covering map $p\colon \widetilde{X}\longrightarrow X$ is called a
\textit{universal covering map} if $\widetilde{X}$ is simply connected.
\end{definition}

\begin{proposition}[Universal covering map of the circle]\label{prop:cubriente-S1}\index{universal covering map of the circle@universal covering map of the circle}
The map $\varepsilon\colon \mathbb R\longrightarrow \mathbb S^1$ defined by $\varepsilon(x)=e^{2\pi i x}$ is a universal covering map of $\mathbb S^1$.
\end{proposition}

\begin{proof}
Let $z_0=e^{2\pi ia}$. If
\[
 U=\left\{e^{2\pi it}\middle|a-\frac{1}{4}<t<a+\frac{1}{4}\right\},
\]
then
\[
 \varepsilon^{-1}(U)
 =\coprod_{k\in\mathbb Z}
 \left(a-\frac{1}{4}+k,a+\frac{1}{4}+k\right),
\]
and each interval projects homeomorphically onto $U$. Thus,
$\varepsilon$ is a covering map. Moreover, $\mathbb R$ is contractible by
the homotopy $H(x,t)=(1-t)x$ and is therefore simply connected.
\end{proof}

\begin{definition}[Evenly covered open subsets]\label{def:regularmente-cubierto}\index{evenly covered open subsets}
Let $p\colon Y\longrightarrow X$ be a covering map of topological spaces.
We say that a connected open subset $U\subseteq X$ is \emph{evenly covered} by $p$ if each connected component $\widetilde{U}$ of $p^{-1}(U)$ projects homeomorphically onto $U$ under the restriction $p\restriction_{\widetilde{U}}\colon \widetilde{U}\longrightarrow U$.
\end{definition}

\begin{theorem}[Arcs covered by the exponential]
\label{teo:epsilon-regular}\index{universal covering map of the circle@universal covering map of the circle!even covering}
Let $U\subsetneq\mathbb S^1$ be an open arc. Then
$\varepsilon^{-1}(U)$ is a disjoint union of open intervals, and the
restriction
\[
 \varepsilon\restriction_{\widetilde U}\colon\widetilde U\longrightarrow U
\]
is a homeomorphism for each component $\widetilde U$ of
$\varepsilon^{-1}(U)$. Thus every proper open arc is evenly
covered in the sense of Definition~\ref{def:regularmente-cubierto}.
\end{theorem}

\begin{proof}
There exist $a\in\mathbb R$ and $0<\ell\leq1$ such that
\[
 U=\{e^{2\pi i t}\mid a<t<a+\ell\}.
\]
Consequently,
\[
 \varepsilon^{-1}(U)
 =\coprod_{k\in\mathbb Z}(a+k,a+\ell+k).
\]
On each interval in this union, $\varepsilon$ is continuous and injective,
and its inverse is a continuous branch of the argument divided by $2\pi$.
Thus each restriction is a homeomorphism onto $U$.
\end{proof}

\section{Complex bundles and stable equivalence}

In what follows, topological vector bundles have finite rank and complex
fibers. Compact base spaces are assumed to be Hausdorff. Noncompact base
spaces are assumed to be locally compact, Hausdorff, and paracompact.
Paracompactness is the hypothesis that allows the construction of the
partitions of unity used to complement and smooth bundles; local
compactness and separation ensure that the one-point compactification is
Hausdorff. The finite-dimensional manifolds considered in this section satisfy these conditions.

\begin{theorem}[Stable complement]\label{teo:complemento-haz-complejo-apendice}
Let $\mathbf{E}\longrightarrow X$ be a complex vector bundle over a compact
Hausdorff space. There exist a bundle $\mathbf{E}'\longrightarrow X$ and an integer $N\geq 0$ such that
\[
\mathbf{E}\oplus \mathbf{E}'\cong X\times\mathbb C^N.
\]
In particular, every stable bundle class can be represented as the
difference between a bundle and a trivial bundle.
\end{theorem}

The proof, including the construction of an embedding of $\mathbf{E}$ into a trivial
bundle by means of a finite partition of unity, can be found in
\cite[Proposition 1.4, pp.~13--14]{HatcherVBKT}. The equivalent formulation in
terms of the Grothendieck group appears in
\cite[Theorem 10.5 and Exercises 10.9--10.12, pp.~262--265]{BleeckerBoossIndex}.

\begin{theorem}[Homotopy invariance of pullback]
\label{teo:homotopia-pullback-apendice}
Let $\mathbf{E}\longrightarrow Y$ be a complex vector bundle and let
$f_0,f_1\colon X\longrightarrow Y$ be homotopic maps. If $X$ is
paracompact, then
\[
f_0^*\mathbf{E}\cong f_1^*\mathbf{E}.
\]
Consequently, the induced homomorphisms
$f_0^*,f_1^*\colon K^0(Y)\longrightarrow K^0(X)$ agree.
\end{theorem}

A proof that trivializes the induced bundle over pieces of $X\times[0,1]$ and
glues the consecutive identifications is given in
\cite[Theorem 1.6, pp.~19--21]{HatcherVBKT}.

\begin{proposition}[Gluing bundles]\label{prop:pegado-haces-apendice}
Let $X=X_0\cup X_1$, with $X_0$ and $X_1$ closed, and let
$\mathbf{E}_i\longrightarrow X_i$ be complex bundles. Each isomorphism
\[
\varphi\colon \mathbf{E}_0\restriction_{X_0\cap X_1}
\longrightarrow \mathbf{E}_1\restriction_{X_0\cap X_1}
\]
determines a bundle $\mathbf{E}_0\cup_\varphi \mathbf{E}_1\longrightarrow X$, unique up to
isomorphism. A homotopy of gluing isomorphisms produces isomorphic
bundles. In particular, bundles over a suspension are described by
transition functions on the equator.
\end{proposition}

The quotient construction and the verification of local triviality are developed
in \cite[Section 1.2, pp.~22--27]{HatcherVBKT}. The parametrized form used
in Bott periodicity can be found in
\cite[Exercise 10.19, pp.~269--270]{BleeckerBoossIndex}.

\section{Relative K-theory and compact support}

For a compact pointed space $(X,x_0)$, define
\[
\widetilde K^0(X)
:=\ker\bigl(i_{x_0}^*\colon K^0(X)\longrightarrow K^0(\{x_0\})\bigr).
\]
If $A\subseteq X$ is closed, set
\[
K^0(X,A):=\widetilde K^0(X/A).
\]
For a locally compact space $Y$, denote its one-point compactification
by $Y^+=Y\cup\{\infty\}$ and set
\[
K_c^0(Y):=\widetilde K^0(Y^+).
\]
If $SX$ is the reduced suspension of a pointed space $X$, we adopt the
convention
\[
\widetilde K^1(X):=\widetilde K^0(SX),
\qquad
K_c^1(Y):=\widetilde K^1(Y^+).
\]

\begin{proposition}[Suspension and product with the line]
\label{prop:k1-producto-r-apendice}
For every locally compact Hausdorff paracompact space $Y$, there is
a natural isomorphism
\begin{equation}
 K_c^1(Y)\cong K_c^0(\mathbb R\times Y).
 \label{eq:k1-producto-r-apendice}
\end{equation}
This isomorphism is compatible with pullbacks by proper maps, the
external product, and extensions associated with open inclusions.
\end{proposition}

\begin{proof}
The one-point compactification of the product admits a pointed homeomorphism
\[
 (\mathbb R\times Y)^+
 \cong S^1\wedge Y^+.
\]
To see this, take the one-point compactification of $\mathbb R$ and consider the
quotient map
\[
 S^1\times Y^+
 \longrightarrow
 \frac{S^1\times Y^+}
 {S^1\times\{\infty\}\,\cup\,\{\infty\}\times Y^+}
 =S^1\wedge Y^+.
\]
Its restriction to $\mathbb R\times Y$ is a homeomorphism onto the
complement of the base point. Neighborhoods of the base point correspond
exactly to complements of compact subsets of $\mathbb R\times Y$; the
universal property of the one-point compactification gives the stated
homeomorphism. Consequently,
\[
 K_c^0(\mathbb R\times Y)
 =\widetilde K^0\bigl((\mathbb R\times Y)^+\bigr)
 \cong\widetilde K^0(S^1\wedge Y^+)
 =\widetilde K^1(Y^+)
 =K_c^1(Y).
\]
All these identifications arise from natural quotient maps. Naturality
with respect to proper maps and compatibility with external products
follow. For an open inclusion, the homeomorphism preserves
the base point and the complement of the open subset, so it also commutes
with extension by zero.
\end{proof}

It is also useful to have a description in which the support remains explicit.
If $K\subseteq Y$ is compact, let
$\mathcal K_K^0(Y)$ denote the group of classes of complexes
\[
\mathbf{E}^0\xrightarrow{\,a\,}\mathbf{E}^1
\]
of complex bundles over $Y$ for which $a_y$ is invertible when
$y\notin K$. Two complexes are identified if, after adding complexes
invertible throughout $Y$, they can be joined by isomorphisms and a homotopy of
bundle homomorphisms that remains invertible outside $K$. For the spaces above,
there is a natural identification
\[
K_c^0(Y)\cong\varinjlim_{K\Subset Y}\mathcal K_K^0(Y).
\]
Support is a notion belonging to this model by complexes, in which the bundle
homomorphism is defined on all of $Y$; no support is assigned to a relative triple
whose bundle homomorphism is defined only on the subspace of the pair.

\begin{theorem}[Exactness, excision, and external product]
\label{teo:exactitud-excision-k-apendice}
The following properties hold.
\begin{enumerate}[label=(\roman*)]
\item For every compact pair $(X,A)$ with $A$ closed, the cyclic sequence
\[
K^0(X,A)\longrightarrow K^0(X)\longrightarrow K^0(A)
\longrightarrow K^1(X,A)\longrightarrow K^1(X)
\longrightarrow K^1(A)\longrightarrow K^0(X,A)
\]
is exact.
\item If $U\subseteq X$ is open and $A\subseteq X$ is closed, with
$\overline U\subseteq\operatorname{int}(A)$, the quotient of pairs induces the
corresponding excision isomorphisms. In particular, if
$j\colon V\hookrightarrow Y$ is an open inclusion of locally compact
spaces, there is an extension homomorphism
\[
j_!\colon K_c^q(V)\longrightarrow K_c^q(Y),
\qquad q\in\{0,1\},
\]
that depends only on the inclusion and agrees with extension by zero for
complexes exact outside a compact set.
More precisely, if $K\Subset V$, restriction to $V$ and extension to
$Y$ induce inverse isomorphisms
\[
\mathcal K_K^0(V)\xrightarrow{\;\cong\;}\mathcal K_K^0(Y).
\]
Taking the limit over compact subsets $K\Subset V$ gives $j_!$.
\item For locally compact spaces $X,Y$, there is a natural associative
external product
\[
\boxtimes\colon K_c^p(X)\otimes K_c^q(Y)
\longrightarrow K_c^{p+q}(X\times Y),
\]
where the degree is taken modulo two. The product commutes with proper
pullbacks and extensions associated with open inclusions.
\end{enumerate}
\end{theorem}

Exactness in degree zero is proved in
\cite[Proposition 2.9, pp.~51--52]{HatcherVBKT}; suspension and
periodicity produce the full cyclic sequence in
\cite[Section 2.2, pp.~51--59]{HatcherVBKT}. The construction using
one-point compactification and the external product can be found in
\cite[Section 10.4 and Theorem 10.16, pp.~266--268]{BleeckerBoossIndex}.
The description by complexes with support and excision in the form just
stated are proved in
\cite[Definition 12.6 and subsequent discussion, pp.~287--289]{BleeckerBoossIndex}.

\begin{theorem}[Relative representation by triples]
\label{teo:modelo-triples-k-apendice}
Let $X$ be a compact Hausdorff space and let $A\subseteq X$ be closed.
The classes of $K^0(X,A)$, with the convention
$K^0(X,\varnothing)=K^0(X)$, are represented by triples $(\mathbf{E},\mathbf{F},\alpha)$, where $\mathbf{E},\mathbf{F}\longrightarrow X$ are complex
bundles and
\[
\alpha\colon \mathbf{E}\restriction_A\longrightarrow \mathbf{F}\restriction_A
\]
is an isomorphism. Two triples represent the same class if, after adding
triples that extend to isomorphisms over $X$, they are related
by an isomorphism and a homotopy of the maps over $A$. Addition is
direct sum, and
\[
[(\mathbf{E},\mathbf{F},\alpha)]+[(\mathbf{F},\mathbf{G},\gamma)]
=[(\mathbf{E},\mathbf{G},\gamma\circ\alpha)].
\]
If $Y=X\setminus A$ and $X/A\cong Y^+$, this model identifies the group with
$K_c^0(Y)$. Under this identification, the model by complexes described above
allows us to choose a bundle homomorphism defined on $Y$ that is invertible outside a
compact set.
\end{theorem}

The construction using the double of $X$, gluing over $A$, and the
difference element can be found in
\cite[Definition 12.6, pp.~287--289]{BleeckerBoossIndex}; the relative model
is developed in \cite[Chapter I, Sections 2--3]{Atiyah1967KTheory}.
For compact Hausdorff spaces, extension of a trivialization from a closed
subset and exactness are established in
\cite[Proposition~2.9 and Section~2.2, pp.~51--55]{HatcherVBKT}.
The proof of Theorem~\ref{teo:modelo-triples-k-relativa} constructs
the isomorphism using the double and verifies the relations,
surjectivity, and injectivity by means of these properties. In particular,
the statement applies to $(X\times\Lambda,A\times\Lambda)$ and its
homotopies for any compact Hausdorff space $\Lambda$, without requiring
$\Lambda$ to be CW.

\begin{theorem}[Smoothing bundles and homomorphisms]
\label{teo:suavizacion-haces-morfismos-apendice}
Let $M$ be a compact smooth manifold with or without boundary.
\begin{enumerate}[label=(\roman*)]
\item Every continuous complex vector bundle over $M$ is isomorphic, as a
continuous bundle, to a smooth complex vector bundle.
\item If $\mathbf{E},\mathbf{F}\longrightarrow M$ are smooth, every continuous bundle homomorphism
$a\colon \mathbf{E}\longrightarrow \mathbf{F}$ can be uniformly approximated by smooth bundle
homomorphisms.
\item If $a$ is an isomorphism and the approximation $a_1$ satisfies
\[
\sup_{x\in M}\lVert a_1(x)-a(x)\rVert
<
\left(\sup_{x\in M}\lVert a(x)^{-1}\rVert\right)^{-1},
\]
then the segment $a_t=(1-t)a+ta_1$ consists of isomorphisms.
In particular, $a$ and $a_1$ define the same relative class.
\end{enumerate}
\end{theorem}

The first assertion and approximation compatible with trivializations are proved
in \cite[p.~101]{Hirsch1976}; see also
\cite[Exercise B.13(b), pp.~719--720]{BleeckerBoossIndex}. The final assertion
follows from
\[
a_t=a\bigl(I+t\,a^{-1}(a_1-a)\bigr)
\]
and the Neumann series, since the norm of the second summand is less than
one uniformly in $x$ and $t$.

\section{Bott periodicity and complex orientation}

Let $\mathbf{H}\longrightarrow\mathbb{CP}^1\cong\mathbb S^2$ be the tautological
line bundle. We identify
$\mathbb R^2$ with $\mathbb S^2\setminus\{\infty\}$ and take as the Bott
class
\[
\beta:=[\mathbf{H}^*]-[\mathbb C]\in
\widetilde K^0(\mathbb S^2)=K_c^0(\mathbb R^2).
\]
With the gluing in Definition~\ref{def:clase-bott}, the tautological bundle is
$\mathbf{H}=\mathbf{E}_1$ and its dual is $\mathbf{H}^*=\mathbf{E}_{-1}$; hence $\beta$ agrees with the class
$b$ used in Chapter~\ref{cap:bott-thom-indice-euclidiano}. This
identification is checked directly, since dualizing a bundle inverts its
gluing function. The description of the tautological bundle through gluing by $z$ is
given in \cite[Exercise B.11, p.~719]{BleeckerBoossIndex}. In the index theory
chapters, this choice is retained in all Thom isomorphisms.
In particular, when a complex line is identified with $\mathbb R^2$ by
$(u,v)\mapsto u+iv$ and $J(u,v)=(-v,u)$ is used, the Thom class satisfies
$\lambda_{\mathbb C}=-\beta$. This convention simultaneously fixes the sign
of the Bott generator, that of the Euclidean model operator, and that of the
normal orientation in the embedding proof.

\begin{theorem}[Complex Bott periodicity]
\label{teo:bott-periodicidad-apendice}
For every locally compact space $X$ and every
$q\in\{0,1\}$, external product with $\beta$ defines a natural isomorphism
\[
\mathcal B_X\colon K_c^q(X)\longrightarrow
K_c^q(X\times\mathbb R^2),
\qquad
\mathcal B_X(a)=a\boxtimes\beta.
\]
In particular,
\[
K_c^0(\mathbb R^{2r})\cong\mathbb Z,
\qquad
K_c^0(\mathbb R^{2r+1})=0,
\]
and $\beta^{\boxtimes r}$ generates the first group.
\end{theorem}

The proof based on families of Wiener--Hopf operators, including the
normalization of the generator and the verification that the two homomorphisms
are inverses, is given in
\cite[Theorems 10.20 and 10.22, pp.~270--272]{BleeckerBoossIndex}. Another complete proof,
formulated using the fundamental product, appears in
\cite[Theorems 2.2 and 2.11, pp.~41--55]{HatcherVBKT}. The original source for
the formulation by complex bundles is
\cite[pp.~229--247]{AtiyahBott1964Periodicity}.

Constructing the pushforward for a nontrivial normal bundle requires
the equivariant version of this normalization. If a compact group $G$ acts
on a space $Y$, let $K_{G,c}^0(Y)$ denote the group defined by
$G$-equivariant complexes that are exact outside a compact set. For a
point, one obtains the complex representation ring
\[
 K_G^0(\{*\})=R(G).
\]

\begin{theorem}[Equivariant Bott normalization]
\label{teo:bott-equivariante-apendice}
Let $G=O(r)$ act on $\mathbb R^r$ by the standard representation, and let
$i\colon\{0\}\hookrightarrow\mathbb R^r$. The equivariant Thom class
\[
 \beta_r:=i_!(1)\in K_{G,c}^0(T^*\mathbb R^r)
\]
has a $G$-equivariant elliptic representative on the compactification
$S^r=\mathbb R^r\cup\{\infty\}$, elementary near $\infty$, whose equivariant
index satisfies
\[
 \operatorname{Ind}_G(\beta_r)=[\mathbb C]
 \quad\text{in }R(G),
\]
where $G$ acts trivially on $\mathbb C$. The result is natural under
formation of bundles associated with a principal $O(r)$-bundle.
\end{theorem}

The construction of the class using the frame bundle, equivariant
descent, and conic cutoffs is developed in
\cite[equations~(12.8)--(12.16), pp.~293--301]{BleeckerBoossIndex}. The
equivariant index identity is formulated in
\cite[p.~295]{BleeckerBoossIndex}. The precise role of
\cite[Proposition~4.4, p.~505]{AtiyahSinger1968I} is as follows: it proves
that axioms \textup{(B1)}, \textup{(B2$'$)}, and
\textup{(B3$''$)} imply the Bott axiom \textup{(B2)}. Within that
proof, the normalization
\[
\operatorname{Ind}_G(i_!1)=1\quad\text{in }R(G),
\]
is established first for blocks of rank one and two, using $O(1)$ and $SO(2)$,
and then in every rank by block decomposition and
multiplicativity. Thus the proposition does not assume the identity
needed here: it obtains that identity as the normalization step recovering
the full Bott axiom.

This normalization agrees with the conventions of the book. In rank one,
we write $z=u+iv$ and $J(u,v)=(-v,u)$; the associated Koszul symbol has
model index $+1$. The gluing function chosen for $b$ has the
opposite orientation, so
$\lambda_{\mathbb C}=-b$, as in
\eqref{eq:thom-linea-menos-bott}. Consequently,
$i_!(1)=\lambda_{\mathbb C}$ has equivariant index equal to the
trivial representation, with no isolated change of sign. This is the
equivariant normalization used for the normal bundle of the embedding.

Now let $\pi\colon \mathbf{V}\longrightarrow X$ be a complex Hermitian bundle of rank
$r$. On $\mathbf{V}$, consider the graded bundles
\[
\pi^*\Lambda^{\mathrm{ev}}\mathbf{V}
:=\bigoplus_{j\ \mathrm{even}}\pi^*\Lambda^j\mathbf{V},
\qquad
\pi^*\Lambda^{\mathrm{odd}}\mathbf{V}
:=\bigoplus_{j\ \mathrm{odd}}\pi^*\Lambda^j\mathbf{V}.
\]
At the point $v\in \mathbf{V}_x$, exterior multiplication by $v$ and its adjoint
define
\[
c_{\mathbf{V}}(v):=\varepsilon(v)-\varepsilon(v)^*\colon
\Lambda^{\mathrm{ev}}\mathbf{V}_x\longrightarrow\Lambda^{\mathrm{odd}}\mathbf{V}_x.
\]
The exterior algebra relations give
\[
c_{\mathbf{V}}(v)^*c_{\mathbf{V}}(v)=\lVert v\rVert^2\operatorname{Id},
\]
so $c_{\mathbf{V}}(v)$ is invertible for $v\neq0$. The triple consisting of the
two bundles and $c_{\mathbf{V}}$ determines the Thom class
\[
\lambda_{\mathbf{V}}\in K_X^0(\mathbf{V}),
\]
where the subscript indicates support on the zero section. If $X$ is compact,
this group is identified with $K_c^0(\mathbf{V})$.
For the trivial line identified with $\mathbb R^2$ by
$(u,v)\mapsto u+iv$ and with $J(u,v)=(-v,u)$, this convention satisfies
\[
\lambda_{\mathbb C}=-\beta.
\]
The verification using transition functions is given in
\eqref{eq:thom-linea-menos-bott}.

\begin{theorem}[Thom isomorphism in K-theory]
\label{teo:thom-k-apendice}
For every complex vector bundle $\pi\colon \mathbf{V}\longrightarrow X$ over a
compact base, the homomorphism
\[
\operatorname{Th}_{\mathbf{V}}\colon K^0(X)\longrightarrow K_c^0(\mathbf{V}),
\qquad
\operatorname{Th}_{\mathbf{V}}(a)=\pi^*a\cdot\lambda_{\mathbf{V}},
\]
is an isomorphism. It is natural with respect to pullbacks, multiplicative
with respect to direct sums, and compatible with the external product. If the
base is merely locally compact, the class with vertical support induces the isomorphism
\[
\operatorname{Th}_{\mathbf{V}}\colon K_c^0(X)\longrightarrow K_c^0(\mathbf{V}).
\]
\end{theorem}

The statement, its normalization, and the comparison with Bott appear in
\cite[Theorem 12.8 and Remark 12.9, pp.~289--290]{BleeckerBoossIndex}.
A proof using Clifford modules is developed in
\cite[Appendix~C, Theorem~C.8 and Remark~C.13]
{LawsonMichelsohn1989}. Independence of
the Hermitian metric also follows directly from the contractibility of
the space of such metrics.

\section{Pushforwards in K-theory}

Let $f\colon X\hookrightarrow Y$ be a proper embedding of smooth manifolds, and
let $\nu_f\longrightarrow X$ be its normal bundle. When $\nu_f$ has a
complex structure, a tubular neighborhood identifies a neighborhood of the
zero section of $\nu_f$ with an open subset $U\subseteq Y$. Define
\[
f_!\colon K_c^0(X)\longrightarrow K_c^0(Y),
\qquad
f_!:=j_!\circ\operatorname{Th}_{\nu_f},
\]
where $j\colon U\hookrightarrow Y$ is the open inclusion.

\begin{theorem}[Properties of the pushforward]
\label{teo:mapa-directo-k-apendice}
The preceding homomorphism is independent of the metric, tubular
neighborhood, and auxiliary identifications. Moreover:
\begin{enumerate}[label=(\roman*)]
\item isotopic embeddings, with compatible complex orientations of their
normal bundles, induce the same map;
\item if $f\colon X\hookrightarrow Y$ and $g\colon Y\hookrightarrow Z$ are
proper embeddings with compatible orientations, then
$(g\circ f)_!=g_!\circ f_!$;
\item for $a\in K_c^0(X)$ and $b\in K_c^0(Y)$, the projection formula
holds:
\[
f_!\bigl(a\cdot f^*b\bigr)=f_!(a)\cdot b;
\]
\item pushforwards commute with external products.
\end{enumerate}
\end{theorem}

The proofs using homotopies of tubular neighborhoods, naturality of
Thom, and excision are given in
\cite[Chapter~III, \S\S12--13]{LawsonMichelsohn1989}. The construction in the
form used for the index theorem appears in
\cite[(12.3)--(12.4), pp.~290--291]{BleeckerBoossIndex}.

For an embedding $f\colon M\hookrightarrow\mathbb R^N$, the tangent embedding
$Tf\colon TM\hookrightarrow T\mathbb R^N$ has normal bundle isomorphic to
\[
\pi^*(\nu_f\oplus\nu_f)\longrightarrow TM.
\]
The rule
\[
J(u,v)=(-v,u)
\]
canonically makes it a complex bundle. Thus, without imposing a
complex structure on $M$, there is a map
\[
(Tf)_!\colon K_c^0(TM)\longrightarrow
K_c^0(T\mathbb R^N)\cong K_c^0(\mathbb R^{2N}).
\]
This is the complex orientation that enters the definition of the general
topological index.

\begin{theorem}[Whitney embedding theorem]\label{teo:whitney-encaje-apendice}
Every Hausdorff second-countable smooth manifold of dimension $n$, with or without boundary, admits
a proper smooth embedding into $\mathbb R^{2n+1}$. If the manifold is compact, its
image is automatically closed.
\end{theorem}

The proof using local embeddings, partitions of unity, and generic
projections can be found in \cite[Theorem 6.15]{LeeS}. The existence and
construction of the required tubular neighborhoods were already proved in
Theorem~\ref{teo:variedades-riemannianas-existencia-de-vecindades-tubulares}.

\section{Characteristic classes and comparison of Thom isomorphisms}

Let $\mathbf{E}\longrightarrow X$ be a complex bundle of rank $r$. The splitting
principle allows us, after a pullback that is injective in cohomology, to work
as though
\[
\mathbf{E}=\mathbf{L}_1\oplus\cdots\oplus \mathbf{L}_r.
\]
If $x_j=c_1(\mathbf{L}_j)$, the elements $x_1,\dots,x_r$ are called formal Chern
roots. The following symmetric expressions descend to well-defined
classes on $X$:
\begin{align*}
c(\mathbf{E})&=\prod_{j=1}^r(1+x_j),\\
\operatorname{ch}(\mathbf{E})&=\sum_{j=1}^r e^{x_j},\\
\operatorname{Td}(\mathbf{E})&=\prod_{j=1}^r\frac{x_j}{1-e^{-x_j}}.
\end{align*}
The last two series are interpreted degree by degree and, in general,
take values in the product $\prod_{j\geq0}H^{2j}(X;\mathbb Q)$.
If $H^k(X;\mathbb Q)=0$ for all sufficiently large $k$, only finitely many
degrees remain, and this product agrees with the direct sum
$H^{\mathrm{ev}}(X;\mathbb Q)$. This holds, in particular, for finite CW
complexes and finite-dimensional smooth manifolds. For a formal
difference, set
\[
\operatorname{ch}([\mathbf{E}]-[\mathbf{F}])
=\operatorname{ch}(\mathbf{E})-\operatorname{ch}(\mathbf{F}).
\]

\begin{theorem}[Chern character]
\label{teo:caracter-chern-apendice}
The Chern character defines a natural transformation of rings
\[
\operatorname{ch}\colon K^0(X)\longrightarrow
\prod_{j\geq0}H^{2j}(X;\mathbb Q).
\]
If $X$ is a finite CW complex, it induces an isomorphism
\[
K^0(X)\otimes_{\mathbb Z}\mathbb Q
\cong H^{\mathrm{ev}}(X;\mathbb Q).
\]
There is also a compactly supported version
\[
\operatorname{ch}_c\colon K_c^0(Y)\longrightarrow
H_c^{\mathrm{ev}}(Y;\mathbb Q)
\]
compatible with external products, proper pullbacks, and extensions
associated with open inclusions.
\end{theorem}

Chern classes, the splitting principle, and the Chern character are
constructed in \cite[Sections 2.3, 3.1 and 4.1, pp.~66--87 and
109--112]{HatcherVBKT}. The compactly supported version and its compatibility
with Thom are developed in
\cite[Section 13.1, pp.~311--315]{BleeckerBoossIndex}.

\begin{theorem}[Cohomological Thom isomorphism and noncompact fundamental class]
\label{teo:thom-cohomologico-borel-moore-apendice}
Let $\pi\colon \mathbf{V}\longrightarrow X$ be an oriented real bundle of rank $r$ over
a compact base of finite CW type. There is a class
$U_{\mathbf{V}}\in H_c^r(\mathbf{V};\mathbb Z)$ whose restriction to each oriented fiber generates
$H_c^r(\mathbb R^r;\mathbb Z)$, and product with this class defines an isomorphism
\[
\operatorname{Th}_{\mathbf{V}}^H\colon H^q(X;R)
\longrightarrow H_c^{q+r}(\mathbf{V};R),
\qquad
\operatorname{Th}_{\mathbf{V}}^H(\omega)=\pi^*\omega\smile U_{\mathbf{V}},
\]
for $R=\mathbb Z,\mathbb Q$ or $\mathbb R$. The Thom class is natural under
oriented pullbacks and multiplicative for oriented direct sums.

If $Y$ is an oriented smooth manifold without boundary of dimension $d$, not necessarily
compact, it has a Borel--Moore fundamental class
$[Y]_{\mathrm{BM}}\in H_d^{\mathrm{BM}}(Y;\mathbb Z)$. The pairing
\[
H_c^d(Y;\mathbb R)\times H_d^{\mathrm{BM}}(Y;\mathbb Z)
\longrightarrow\mathbb R
\]
agrees, under the compactly supported de Rham theorem, with
$\displaystyle [\omega]\mapsto\int_Y\omega$.
\end{theorem}

The construction of the cohomological orientation, the compactly supported
Thom isomorphism, its naturality, and evaluation on the fundamental class
are developed in
\cite[Section 13.1, pp.~311--315]{BleeckerBoossIndex}.

With the Thom class fixed above, comparison of the Thom isomorphisms in
K-theory and cohomology introduces the Todd factor. Under the
exterior multiplication convention of $\lambda_{\mathbf{V}}$, for a complex bundle
of rank $r$ one has
\[
\operatorname{ch}_c(\lambda_{\mathbf{V}})
=(-1)^r\operatorname{Th}^{H}_{\mathbf{V}}
\bigl(\operatorname{Td}(\overline{\mathbf{V}})^{-1}\bigr),
\]
where $\operatorname{Th}^{H}_{\mathbf{V}}$ is the Thom isomorphism in cohomology
with rational coefficients.

\begin{theorem}[Cohomological formula for the topological index]
\label{teo:formula-cohomologica-indice-apendice}
Let $M$ be a closed smooth manifold of dimension $n$, and let
$u\in K_c^0(T^*M)$. Identify $T^*M$ with $TM$ using a metric and
orient $TM$ by the almost complex structure that takes the horizontal
directions as real and the vertical directions as imaginary. If
$a=(\sharp_g)_*u\in K_c^0(TM)$, then
\[
\operatorname{ind}_{\mathrm t}(u)
=(-1)^n
\int_{TM}
\operatorname{ch}_c(a)\smile
\pi^*\operatorname{Td}(TM\otimes_{\mathbb R}\mathbb C).
\]
The integrand has compact support in the fiber direction; the
integral means evaluation of its degree-$2n$ component on the
fundamental class of $TM$.
\end{theorem}

The full comparison of the two Thom classes, including the conjugate bundle
and the sign
$(-1)^n$ for these conventions, is given in
\cite[Theorem 13.1, pp.~313--315]{BleeckerBoossIndex}. This is the formulation
used in the chapter devoted to the Atiyah--Singer theorem.

\section{De Rham cohomology, spin structures, and the Euler class}

The results in this section are used only to identify the indices of
geometric operators topologically. Their proofs belong to algebraic and
differential topology and are not reproduced here.

We first recall the obstruction that arises for spin structures.
For a real vector bundle $\mathbf{V}\to X$ of rank $r$ over a paracompact
space, its \emph{Stiefel--Whitney classes} are the classes
\[
w_j(\mathbf{V})\in H^j(X;\mathbb Z_2),
\qquad 0\leq j\leq r,
\]
characterized by naturality, $w_0(\mathbf{V})=1$, the Whitney formula
\[
w(\mathbf{V}\oplus \mathbf{W})=w(\mathbf{V})\smile w(\mathbf{W}),
\qquad
w(\mathbf{V}):=1+w_1(\mathbf{V})+\cdots+w_r(\mathbf{V}),
\]
and the normalization $w_1(\boldsymbol{\gamma}^1)\neq0$ for the real tautological bundle
$\boldsymbol{\gamma}^1\to\mathbb{RP}^{\infty}$. In particular,
\[
w_2(\mathbf{V})\in H^2(X;\mathbb Z_2)
\]
is the degree-two component of the total class. The vanishing of $w_1(\mathbf{V})$
is equivalent to orientability of $\mathbf{V}$; after choosing an orientation
and a continuous Euclidean metric on the fibers, $w_2(\mathbf{V})$ is the
obstruction to lifting the principal bundle of oriented orthonormal frames
through the covering map
$\operatorname{Spin}(r)\to\operatorname{SO}(r)$.

The universal construction of these classes, the orientation criterion, and the
interpretation of $w_2$ as the spin obstruction can be found in
\cite[Definition~17.18 and Proposition~17.20, pp.~527--529]{BleeckerBoossIndex}.
We will use only those topological results; Clifford algebras, the
spin covering map, and spinor bundles are constructed geometrically in
Section~\ref{sec:clifford-spin-dirac}.

\begin{theorem}[De Rham theorem and Euler characteristic]
\label{teo:de-rham-euler-apendice}
Let $M$ be a smooth manifold with or without boundary. Integration of forms over smooth chains
induces, for each $q$, a natural isomorphism
\[
H_{\mathrm{dR}}^q(M;\mathbb C)
\cong H_{\mathrm{sing}}^q(M;\mathbb C).
\]
If $M$ is a closed manifold of dimension $n$, these spaces are
finite-dimensional and
\[
\chi(M)=\sum_{q=0}^{n}(-1)^q
\dim_{\mathbb C}H_{\mathrm{dR}}^q(M;\mathbb C).
\]
\end{theorem}

The identification with singular cohomology and its relation to the
Euler characteristic are presented in
\cite[Section 13.4, in particular (13.18), pp.~319--323]{BleeckerBoossIndex};
the Hodge--de Rham decomposition that complements this identification appears
in \cite[Theorem 17.62 and Corollary 17.63, pp.~602--604]{BleeckerBoossIndex}.

\begin{theorem}[Euler class and Euler characteristic]
\label{teo:euler-caracteristica-apendice}
Let $\mathbf{V}\longrightarrow X$ be an oriented real vector bundle of rank $r$ over
a compact base. If
$\mathbf{s}_0\colon X\longrightarrow \mathbf{V}$ is the zero section and $U_{\mathbf{V}}$ is the cohomological
Thom class, the Euler class is defined by
\[
 e(\mathbf{V}):=\mathbf{s}_0^*U_{\mathbf{V}}\in H^r(X;\mathbb Z).
\]
This class is natural under oriented pullbacks and multiplicative with respect
to oriented direct sums. If $M$ is a closed oriented smooth manifold of
dimension $n$, then
\[
 \bigl\langle e(TM),[M]\bigr\rangle=\chi(M).
\]
\end{theorem}

The construction using the zero section and evaluation on the fundamental
class can be found in
\cite[Theorem~13.6(d), pp.~320--321]{BleeckerBoossIndex} when the dimension
is even. In odd dimensions, Poincaré duality pairs Betti numbers in
complementary degrees and gives $\chi(M)=0$. On the other hand,
$-\operatorname{Id}_{TM}$ reverses the orientation of each fiber, so
naturality of the Thom class implies $e(TM)=-e(TM)$. Since
$H^n(M;\mathbb Z)$ is free for a closed oriented manifold, it follows
that $e(TM)=0$. This proves the statement in every dimension. This equality
is the general topological form of the Gauss--Bonnet--Chern theorem.

Let
\[
0\longrightarrow\mathbb Z
\xrightarrow{\ \times2\ }\mathbb Z
\longrightarrow\mathbb Z_2
\longrightarrow0
\]
be the exact sequence of coefficients, and let
\[
\beta\colon H^2(X;\mathbb Z_2)\longrightarrow H^3(X;\mathbb Z)
\]
be its Bockstein homomorphism. For an oriented real bundle $\mathbf{V}\to X$,
define the third integral Stiefel--Whitney class by
\[
W_3(\mathbf{V}):=\beta\bigl(w_2(\mathbf{V})\bigr).
\]
By exactness, $W_3(\mathbf{V})=0$ if and only if $w_2(\mathbf{V})$ is the reduction modulo two
of a class in $H^2(X;\mathbb Z)$.

\begin{theorem}[$\operatorname{Spin}^c$ orientation]
\label{teo:spinc-apendice}
An oriented real bundle $\mathbf{V}\longrightarrow X$ over a paracompact base
admits a $\operatorname{Spin}^c$ structure if and only if
$W_3(\mathbf{V})=0$; equivalently, if and only if $w_2(\mathbf{V})$ is the reduction modulo
two of an integral class of degree two. The determinant line $\mathbf{L}$ of a
chosen structure satisfies
\[
c_1(\mathbf{L})\bmod2=w_2(\mathbf{V}).
\]
If the rank of $\mathbf{V}$ is even, Clifford multiplication between the
half-spinor bundles defines a Thom class in $K_X^0(\mathbf{V})$, the group with
support on the zero section; when $X$ is compact, this is identified with
$K_c^0(\mathbf{V})$. Every complex bundle
$\mathbf{E}$, regarded as an oriented real bundle, has a canonical
$\operatorname{Spin}^c$ structure induced by
$U(r)\to\operatorname{Spin}^c(2r)$; its determinant line is
$\det_{\mathbb C}\mathbf{E}$ and its Thom class agrees with the class constructed
using the exterior algebra.
\end{theorem}

The existence criterion is proved in
\cite[Theorem 18.35, pp.~655--656]{BleeckerBoossIndex} and in
\cite[Appendix~D, Theorem~D.2]{LawsonMichelsohn1989}. The canonical
orientation of a complex bundle and comparison of the two Thom classes are
developed in
\cite[Appendix~D, Example~D.6; Appendix~C,
Theorems~C.8 and C.12]{LawsonMichelsohn1989}.

\begin{theorem}[Spin obstruction and half-spinor bundles]
\label{teo:spin-semiespinores-apendice}
Let $\mathbf{V}\to X$ be an oriented real bundle of rank $n$ over a
paracompact base. After choosing a continuous Euclidean metric on its
fibers, its oriented orthonormal frame bundle admits a principal lift
with group $\operatorname{Spin}(n)$ if and only if
\[
w_2(\mathbf{V})=0\in H^2(X;\mathbb Z_2).
\]
When such a lift exists, the set of equivalence classes of spin structures is
a principal homogeneous space for $H^1(X;\mathbb Z_2)$. If $n=2m$, each
spin structure determines Hermitian half-spinor bundles
$\mathbf{S}^+$ and $\mathbf{S}^-$, and Clifford multiplication by a nonzero vector
defines an isomorphism
\[
c(v)\colon \mathbf{S}_x^+\longrightarrow \mathbf{S}_x^-,
\qquad v\in \mathbf{V}_x\setminus\{0\}.
\]
The class of this isomorphism is the spinorial Thom orientation of
$\mathbf{V}$ in complex $K$-theory, with support on the zero section. If the formal
roots of
$\mathbf{V}\otimes_{\mathbb R}\mathbb C$ are ordered as
$x_1,-x_1,\ldots,x_m,-x_m$, then
\[
\operatorname{ch}(\mathbf{S}^+)-\operatorname{ch}(\mathbf{S}^-)
=\prod_{j=1}^m
\left(e^{-\frac{x_j}{2}}-e^{\frac{x_j}{2}}\right),
\qquad
\widehat A(\mathbf{V})=\prod_{j=1}^m
\frac{\frac{x_j}{2}}{\sinh\!\left(\frac{x_j}{2}\right)}.
\]
Here the orientation is chosen so that
$e(\mathbf{V})=x_1\cdots x_m$; the order of the two exponents is determined by
the chirality convention in
Section~\ref{sec:clifford-spin-dirac}.
\end{theorem}

The construction of $w_2$ and the existence criterion can be found in
\cite[Definition 17.18 and Proposition 17.20, pp.~527--529]{BleeckerBoossIndex}.
The bundles $\mathbf{S}^\pm$, Clifford multiplication, and the twisted Dirac operator
are constructed in
\cite[Definitions 17.22--17.23, pp.~531--532]{BleeckerBoossIndex}; the
comparison of its symbol class with the Thom class is developed in
\cite[Appendix~C, Theorem~C.12]{LawsonMichelsohn1989}.
The calculation of the weights of the half-spinor representations that produces
$\operatorname{ch}(\mathbf{S}^+)-\operatorname{ch}(\mathbf{S}^-)$ is given in
\cite[pp.~519--522]{BleeckerBoossIndex}; the class $\widehat A$ and its relation
to the Dirac index appear in
\cite[equations (17.43)--(17.45), pp.~545--546]{BleeckerBoossIndex}.

\begin{theorem}[Chern defect of the spinorial Thom class]
\label{teo:defecto-chern-thom-spinorial-apendice}
Let $\mathbf{V}\to X$ be a real spin bundle of rank $2m$ over a compact base, and let
\[
s(\mathbf{V})=[\pi^*\mathbf{S}^+,\pi^*\mathbf{S}^-;c(v)]\in K_c^0(\mathbf{V})
\]
be its spinorial Thom class with the preceding orientations. Then
\begin{equation}
(\operatorname{Th}_{\mathbf{V}}^H)^{-1}
\operatorname{ch}_c\bigl(s(\mathbf{V})\bigr)
=(-1)^m\widehat A(\mathbf{V})^{-1}.
\label{eq:defecto-chern-thom-spinorial}
\end{equation}
\end{theorem}

The formula, including the sign determined by the choice of chirality,
is Proposition~12.5 in Chapter~III, \S12, of
\cite{LawsonMichelsohn1989}.

\begin{proposition}[Topological facts for extinction by surgery]
\label{prop:topologia-extincion-cirugia}
Let $M_1,\ldots,M_k$ be closed connected manifolds of dimension at least
three. Then
\[
\pi_1(M_1\#\cdots\#M_k)
\cong
\pi_1(M_1)*\cdots *\pi_1(M_k).
\]
If $N=S^3/\Gamma$ is a three-dimensional spherical space form, where
$\Gamma$ is finite and acts freely, then
$\pi_1(N)\cong\Gamma$; in particular, $N$ is simply connected if and only
if $N\cong S^3$. If $Z\to S^1$ is a bundle with fiber $S^2$, trivial or
twisted, then $\pi_1(Z)\cong\mathbb Z$.

Consequently, if a closed connected manifold $M$ is a connected sum of
spherical space forms and $S^2$-bundles over $S^1$, and
$\pi_1(M)=0$, there are no bundle factors and all the groups
$\Gamma$ are trivial. Since $S^3$ is the identity for connected sum,
$M\cong S^3$.
\end{proposition}

\begin{proof}
The first identity is Seifert--van Kampen applied to the neck of the
connected sum. The second follows because $S^3$ is the universal covering space of
$S^3/\Gamma$, and the third follows from the homotopy exact sequence of the bundle,
since $\pi_1(S^2)=0$. A free product is trivial only if each of its
factors is trivial; this gives the conclusion.
\end{proof}

Every simply connected manifold is orientable. Moreover, a closed two-sided
surface in an orientable $3$-manifold inherits an orientation from
the ambient orientation and a global normal; hence such a manifold
contains no two-sided embedded $\mathbb{RP}^2$. The topological effect
of surgeries and the reconstruction by connected sums that we will use are
formulated in \cite[Chapter~15]{morgan2007ricci}; the extinction used at
the end is treated in Chapters~18--19 of the same reference.

\begin{theorem}[Gauss--Bonnet for surfaces]
\label{teo:gauss-bonnet-superficies-apendice}
Let $(M,\mathbf{g})$ be a closed Riemannian surface and let $K_{\mathbf{g}}$ be its Gaussian
curvature. Then
\[
\int_M K_{\mathbf{g}}\,d\lambda_{\mathbf{g}}=2\pi\chi(M).
\]
Since $R_{\mathbf{g}}=2K_{\mathbf{g}}$, it follows that
\[
\int_M R_{\mathbf{g}}\,d\lambda_{\mathbf{g}}=4\pi\chi(M).
\]
\end{theorem}

The formulation using the Euler class can be found in
\cite[Theorem 13.6(d), pp.~320--321]{BleeckerBoossIndex}; a proof
using the Gauss--Bonnet--Chern form appears in
\cite[Theorem 17.68, pp.~612--613]{BleeckerBoossIndex}.
If $M$ is nonorientable, apply the formula to its orientable
two-sheeted covering space:
the integral and the Euler characteristic are both multiplied by two, so
the same identity descends to $M$.

\chapter{Mathematical analysis}\label{ap:analisis-matematico}

This appendix collects the results from real analysis, metric spaces, functional analysis, measure theory, and Lebesgue spaces used throughout the book. Its purpose is not to replace a systematic treatment of these subjects, but to fix notation, recall frequently used tools, and provide precise internal references for the main chapters. More extensive treatments can be found in \cite{Amann,Grabinsky,Bachman,Monica,DrabekMilota2013,brezis_functional_2011,kesavan_functional_2023,Rudin1991,reed_simon_vol1}.

The order follows the natural dependencies between the topics. We begin with some facts from real analysis and metric spaces; continue with normed spaces, linear operators, and Hilbert spaces; and then study measure theory, the Lebesgue and Bochner integrals, and $L^{p}(X)$ spaces. Finally, we introduce topological vector spaces and weak topologies, which play an essential role in the compactness and semicontinuity arguments used in the calculus of variations.

\section{Real analysis and metric spaces}\label{ap:analisis-real-espacios-metricos}

We begin with some basic facts from real analysis and metric spaces. These results will be used mainly to control local arguments, completions, and convergence.
\begin{proposition}\label{derivada y límite}
 Let $f\colon (a,b)\longrightarrow \mathbb{R}$ be a function differentiable on $(a,b)\setminus \{c\}$, with $c\in (a,b)$, and continuous at $(a,b)$, such that $\displaystyle\lim_{x\to c^{+}}f'(x)=L=\displaystyle\lim_{x\to c^{-}}f'(x)$. Then $f$ is differentiable at $c$ and $f'(c)=L$.
\end{proposition}
\begin{proof}
Let $x\in(a,b)\setminus\{c\}$. If $x>c$, the mean value theorem applied to $f$ on $[c,x]$ gives a point $\xi_x\in(c,x)$ such that
\[
\frac{f(x)-f(c)}{x-c}=f'(\xi_x).
\]
As $x\to c^+$, necessarily $\xi_x\to c^+$, and therefore
\[
\lim_{x\to c^+}\frac{f(x)-f(c)}{x-c}=L.
\]
If $x<c$, the same argument applied to the interval $[x,c]$ produces $\xi_x\in(x,c)$ and gives
\[
\frac{f(x)-f(c)}{x-c}=f'(\xi_x).
\]
Now $\xi_x\to c^-$ as $x\to c^-$, so the left-hand limit of the difference quotient is also $L$. The two one-sided limits agree; consequently, $f$ is differentiable at $c$ and $f'(c)=L$.
\end{proof}

\begin{proposition}\label{normas en R^{n}}
 For every $x\in\mathbb{R}^{n}$, the following hold:
 \begin{enumerate}[label=(\alph*)]
 \item $\|x\|_{r}\leq\|x\|_{s}$ if $1\leq s\leq r\leq \infty$,
 \item $\|x\|_{s}\leq n^{\displaystyle\frac{r-s}{sr}}\|x\|_{r}$ if $1\leq s\leq r<\infty$,
 \item $\|x\|_{\infty}\leq \|x\|_{s}\leq n^{\displaystyle\frac{1}{s}}\|x\|_{\infty}$ if $1\leq s < \infty$.
 \end{enumerate}
 In particular, all norms $\|\cdot\|_{p}$ with $p\in [1,\infty]$ are equivalent on $\mathbb{R}^{n}$.
\end{proposition}
\begin{proof}
If $x=0$, all the inequalities are immediate. Suppose that $x\neq0$. First let $1\leq s\leq r<\infty$ and set
\[
a_i:=\frac{|x_i|}{\|x\|_s},\qquad 1\leq i\leq n.
\]
Then $0\leq a_i\leq1$ and $\displaystyle \sum_{i=1}^n a_i^s=1$. Since $r\geq s$, we have $a_i^r\leq a_i^s$, and hence
\[
\frac{\|x\|_r^r}{\|x\|_s^r}
=\sum_{i=1}^n a_i^r
\leq\sum_{i=1}^n a_i^s
=1.
\]
This proves $\|x\|_r\leq\|x\|_s$. For $r=\infty$, the same inequality follows directly from $\displaystyle |x_i|^s\leq\displaystyle\sum_{j=1}^n|x_j|^s$ for each $i$.

For the reverse inequality, the case $s=r$ is trivial. If $1\leq s<r<\infty$, apply Hölder's inequality with conjugate exponents $\displaystyle\frac{r}{s}$ and $\displaystyle\frac{r}{r-s}$:
\[
\sum_{i=1}^n|x_i|^s
\leq
\left(\sum_{i=1}^n|x_i|^r\right)^{\frac{s}{r}}
\left(\sum_{i=1}^n1^{\,\frac{r}{r-s}}\right)^{\frac{r-s}{r}}
=\|x\|_r^s n^{1-\frac{s}{r}}.
\]
Taking $s$th roots gives
\[
\|x\|_s\leq n^{\frac{1}{s}-\frac{1}{r}}\|x\|_r
=n^{\frac{r-s}{sr}}\|x\|_r.
\]
Finally, for $1\leq s<\infty$,
\[
\|x\|_\infty^s
\leq\sum_{i=1}^n|x_i|^s
\leq n\|x\|_\infty^s,
\]
and taking roots gives part (c). The inequalities in parts (a)--(c) show that the identity between any two of these normed spaces is continuous in both directions; hence all the indicated norms are equivalent.
\end{proof}

\begin{proposition}\label{restriction of C1 function to compact is lipschitz}
 If $\Omega\subseteq\mathbb{R}^{n}$ is open and $f\colon \Omega\longrightarrow \mathbb{R}^{m}$ is of class $C^{1}$, then for every compact convex set $K\subseteq \Omega$, the restriction $f\restriction_{K}$ is Lipschitz.
\end{proposition}
\begin{proof}
If $K=\varnothing$, there is nothing to prove. Since $Df$ is continuous and $K$ is compact, there exists
\[
L:=\max_{z\in K}\|Df(z)\|_{\mathrm{op}}<\infty.
\]
Given $x,y\in K$, convexity of $K$ implies that the segment
$\gamma(t):=x+t(y-x)$, $0\leq t\leq1$, is contained in $K$. The chain rule and the fundamental theorem of calculus, applied componentwise to $f\circ\gamma$, give
\[
f(y)-f(x)
=\int_0^1\frac{d}{dt}f(\gamma(t))\,dt
=\int_0^1Df(\gamma(t))(y-x)\,dt.
\]
Consequently,
\[
\|f(y)-f(x)\|
\leq\int_0^1\|Df(\gamma(t))\|_{\mathrm{op}}\|y-x\|\,dt
\leq L\|y-x\|.
\]
Thus $f\restriction_K$ is $L$-Lipschitz.
\end{proof}
\begin{remark}\label{obs:b1-analisis-real-y-metricos-hipotesis-convexidad-reemplazarse-existencia-control-uni}
 The convexity hypothesis can be replaced by uniform control of the lengths of curves in $K$ joining any two points of $K$. In particular, the result applies locally on closed balls contained in $\Omega$, which is the case used most frequently.
\end{remark}
\begin{definition}\label{def:b1-analisis-real-y-metricos-metrica}\index{metric@metric}
 Let $X$ be a set. A \textbf{metric} on $X$
 is a function $d\colon X\times X\longrightarrow [0,\infty)$ such that:
\begin{enumerate}[label=(\alph*)]
 \item $d(x,y)=0$ if and only if $x=y$.
 \item $d(x,y)=d(y,x)$ for all $x,y\in X$.
 \item $d(x,z)\leq d(x,y)+d(y,z)$ for all $x,y,z\in X$.
\end{enumerate} A \textbf{metric space} is a pair $(X,d)$, where $X$ is a set and $d$ is a metric on $X$.
\end{definition}

\begin{theorem}[McShane extension theorem]\label{teo:extension-mcshane}\index{McShane extension@McShane extension}\index{Lipschitz function@Lipschitz function!extension}
Let $(X,d)$ be a metric space, $A\subseteq X$ a nonempty subset, and $f\colon A\longrightarrow\mathbb R$ an $L$-Lipschitz function, with $L\geq0$. Then there is an extension $F\colon X\longrightarrow\mathbb R$ of $f$ that is also $L$-Lipschitz. An explicit choice is
\[
F(x):=\inf_{a\in A}\bigl(f(a)+L\,d(x,a)\bigr).
\]
\end{theorem}

\begin{proof}
Fix $a_0\in A$. For any $x\in X$ and $a\in A$, the Lipschitz property and the triangle inequality give
\[
f(a)+L\,d(x,a)
\geq
f(a_0)-L\,d(a,a_0)+L\,d(x,a)
\geq
f(a_0)-L\,d(x,a_0).
\]
On the other hand, taking $a=a_0$ gives $F(x)\leq f(a_0)+L\,d(x,a_0)$. Thus the infimum defining $F(x)$ is finite.

If $x\in A$, the inequality $f(x)\leq f(a)+L\,d(x,a)$ holds for every $a\in A$, so $f(x)\leq F(x)$; choosing $a=x$ gives the reverse inequality. Thus $F\restriction_A=f$. Finally, for any $x,y\in X$ and $a\in A$,
\[
f(a)+L\,d(x,a)
\leq
f(a)+L\,d(y,a)+L\,d(x,y).
\]
Taking the infimum over $a$ gives $F(x)\leq F(y)+L\,d(x,y)$. Interchanging $x$ and $y$ yields $|F(x)-F(y)|\leq L\,d(x,y)$.
\end{proof}

\begin{definition}\label{def:b1-analisis-real-y-metricos-esimo-termino}\index{sequence@sequence}
 A sequence in a metric space $X=(X,d)$ is a function $x\colon \mathbb{N}\longrightarrow X$. The value of this function at $k$ is called the \textbf{$k$th term} of the sequence and is denoted by $x_{k}$. A sequence is usually denoted by $(x_{k})$ or $\{x_{k}\}_{k\in\mathbb{N}}$, or by similar notation according to context.

 We say that $(x_{k})$ converges to a point $x\in X$ if for every $\varepsilon>0$ there exists $N\in\mathbb{N}$ such that if $k\geq N$, then $d(x_{k},x)<\varepsilon$.

 We say that a sequence $(x_{k})$ is Cauchy if for every $\varepsilon>0$ there exists $N\in\mathbb{N}$ such that if $k,l\geq N$, then $d(x_{k},x_{l})<\varepsilon$.
\end{definition}
In every metric space, convergent sequences are Cauchy, but a Cauchy sequence need not converge. This leads to the following definition:
\begin{definition}\label{def:b1-analisis-real-y-metricos-espacio-metrico-completo}\index{complete metric space@complete metric space}
 A metric space $(X,d)$ is a \textbf{complete metric space} if every Cauchy sequence in $X$ converges in $X$.
\end{definition}
\begin{definition}\label{def:b1-analisis-real-y-metricos-isometria}\index{isometry@isometry}
 Let $(X,d_{X})$ and $(Y,d_{Y})$ be metric spaces. A function $\phi\colon X\longrightarrow Y$ is an \textbf{isometry} if $d_{Y}(\phi(x),\phi(y))=d_{X}(x,y)$ for all $x,y\in X$.
\end{definition}
\begin{definition}\label{def:b1-analisis-real-y-metricos-espacio-metrico-espacio-metrico-completo-completacion}\index{completion@completion}
 Let $X$ be a metric space. A complete metric space $\overline{X}$ is called a \textbf{\textit{completion}} of $X$ if there is an isometry $i\colon X\longrightarrow \overline{X}$ such that $i(X)$ is dense in $\overline{X}$.
\end{definition}
We will prove that every metric space admits a completion up to isometry, but first we need some lemmas.
\begin{lemma}\label{lema para completacion}
 Let $(X,d)$ be a metric space, $(b_{n})$ a Cauchy sequence in $X$, and $(a_{n})$ a sequence in $X$ such that $d(a_{n},b_{n})<\displaystyle\frac{1}{n}$ for every $n\in\mathbb{N}$. Then:
 \begin{enumerate}[label=(\alph*)]
 \item $(a_{n})$ is a Cauchy sequence in $X$.
 \item $a_{n}\to x\in X$ if and only if $b_{n}\to x\in X$.
 \end{enumerate}
\end{lemma}
\begin{proof}
 For $(a)$, the triangle inequality gives \[d(a_{n},a_{m})\leq d(a_{n},b_{n})+d(b_{n},b_{m})+d(b_{m},a_{m})\], and by hypothesis there exists $N_{1}\in\mathbb{N}$ such that if $n,m\geq N_{1}$, then $d(a_{n},b_{n})<\displaystyle\frac{\varepsilon}{3}$ and $d(b_{m},a_{m})<\displaystyle\frac{\varepsilon}{3}$.

 Finally, since $(b_{n})$ is Cauchy, there exists $N_{2}\in\mathbb{N
 }$ such that if $n,m\geq N_{2}$, then $d(b_{n},b_{m})<\displaystyle\frac{\varepsilon}{3}$. Thus, if $n,m\geq\displaystyle\max\{N_{1},N_{2}\}$, then $d(a_{n},a_{m})<3\displaystyle\frac{\varepsilon}{3}=\varepsilon$.

 For $(b)$, if $a_{n}\to x$, then $d(a_{n},x)\to 0$ whenever $n\to\infty$. By the triangle inequality, $d(b_{n},x)\leq d(b_{n},a_{n})+d(a_{n},x)<\displaystyle\frac{1}{n}+d(a_{n},x)$. Taking $n\to\infty$, we conclude that $d(b_{n},x)\to 0$, and hence $b_{n}\to x$. The converse is proved analogously.
\end{proof}
\begin{theorem}\label{completación de espacio métrico}\index{completion of a metric space@completion of a metric space}
 Every metric space admits a completion, unique up to isometry.
\end{theorem}
\begin{proof}
Let $CS(X)$ be the set of Cauchy sequences in $X$. On $CS(X)$, define the relation
\[
(x_n)\sim(y_n)
\quad\Longleftrightarrow\quad
\lim_{n\to\infty}d(x_n,y_n)=0.
\]
Reflexivity and symmetry are immediate. For transitivity, if $(x_n)\sim(y_n)$ and $(y_n)\sim(z_n)$, the triangle inequality gives $0\leq d(x_n,z_n)\leq d(x_n,y_n)+d(y_n,z_n)$; passing to the limit yields $(x_n)\sim(z_n)$. Thus $\sim$ is an equivalence relation.

Let $\overline X:=CS(X)/\sim$. For $[(x_n)],[(y_n)]\in\overline X$, set
\[
\overline d\bigl([(x_n)],[(y_n)]\bigr)
:=
\lim_{n\to\infty}d(x_n,y_n).
\]
The limit exists, since for any $m,n\in\mathbb N$ the triangle inequality implies
\[
\bigl|d(x_n,y_n)-d(x_m,y_m)\bigr|
\leq d(x_n,x_m)+d(y_n,y_m),
\]
and the right-hand side tends to zero because $(x_n)$ and $(y_n)$ are Cauchy. This estimate shows that $(d(x_n,y_n))$ is a Cauchy sequence in $\mathbb R$.

The definition is independent of the representatives. If $(x_n)\sim(x_n')$ and $(y_n)\sim(y_n')$, then
\[
\bigl|d(x_n,y_n)-d(x_n',y_n')\bigr|
\leq d(x_n,x_n')+d(y_n,y_n'),
\]
so the two limits agree. The metric properties follow directly from the corresponding properties of $d$; in particular, the triangle inequality is obtained by passing to the limit in $d(x_n,z_n)\leq d(x_n,y_n)+d(y_n,z_n)$. Thus $(\overline X,\overline d)$ is a metric space.

For each $x\in X$, let $x^*:=[(x,x,\ldots)]$ and define $i\colon X\longrightarrow\overline X$ by $i(x):=x^*$. If $x,y\in X$, then
\[
\overline d(i(x),i(y))
=
\lim_{n\to\infty}d(x,y)
=
d(x,y),
\]
so $i$ is an isometric embedding. Its image is dense. Given a class $[(x_n)]\in\overline X$, the sequence $(x_m^*)$ satisfies
\[
\overline d(x_m^*,[(x_n)])
=
\lim_{n\to\infty}d(x_m,x_n)
\longrightarrow 0
\]
as $m\to\infty$, since $(x_n)$ is Cauchy.

We now prove that $\overline X$ is complete. Let $(a_n)$ be a Cauchy sequence in $\overline X$. By density of $i(X)$, for each $n$ we can choose $x_n\in X$ such that $\overline d(i(x_n),a_n)<\displaystyle\frac{1}{n}$. Lemma~\ref{lema para completacion} shows that $(i(x_n))$ is Cauchy in $\overline X$. Since $i$ is an isometry, $(x_n)$ is Cauchy in $X$, and the construction of $\overline X$ gives $i(x_n)\to[(x_n)]$. Another application of the same lemma yields $a_n\to[(x_n)]$. This proves that $\overline X$ is a completion of $X$.

It remains to verify uniqueness. Let $(Y,d_Y)$ be another completion of $X$, and let $j\colon X\longrightarrow Y$ be an isometric embedding with dense image. For each $y\in Y$, choose a sequence $(x_n)$ in $X$ such that $j(x_n)\to y$ and define
\[
G(y):=[(x_n)]\in\overline X.
\]
The map $G$ is well defined: if also $j(x_n')\to y$, then $d(x_n,x_n')=d_Y(j(x_n),j(x_n'))\to0$, so both sequences determine the same class. Moreover, for $y,y'\in Y$ and sequences $j(x_n)\to y$, $j(x_n')\to y'$, we have
\[
\overline d(G(y),G(y'))
=
\lim_{n\to\infty}d(x_n,x_n')
=
\lim_{n\to\infty}d_Y(j(x_n),j(x_n'))
=
d_Y(y,y').
\]
Therefore $G$ is an isometry. It is surjective because every class $[(x_n)]\in\overline X$ is obtained from the limit in $Y$ of the Cauchy sequence $(j(x_n))$. Consequently, $G$ is a bijective isometry and the completion is unique up to isometry.
\end{proof}

Completeness yields one of the structural principles underlying much of functional analysis. We first recall the corresponding topological notion.

\begin{definition}\label{def:espacio-de-baire}\index{Baire space}
A topological space is called a \textbf{Baire space} if every countable intersection of open dense subsets is dense. Equivalently, no nonempty open subset can be written as a countable union of closed sets with empty interior.
\end{definition}

\begin{theorem}[Baire category theorem]\label{teo:categoria-de-baire}\index{Baire category theorem@Baire category theorem}
Every complete metric space is a Baire space.
\end{theorem}

\begin{proof}
Let $(X,d)$ be a complete metric space, and let $(G_n)_{n\in\mathbb N}$ be open dense subsets of $X$. We must prove that $\displaystyle\bigcap_{n=1}^{\infty}G_n$ meets every nonempty open subset $U\subseteq X$. Choose $x_1\in U\cap G_1$ and $r_1\in(0,1)$ such that $\overline{B}_d(x_1,r_1)\subseteq U\cap G_1$. Once $x_n$ and $r_n$ have been chosen, density of $G_{n+1}$ allows us to choose $x_{n+1}\in B_d(x_n,\displaystyle\frac{r_n}{2})\cap G_{n+1}$ and $r_{n+1}\in(0,2^{-n-1})$ so that $\overline{B}_d(x_{n+1},r_{n+1})\subseteq B_d(x_n,\displaystyle\frac{r_n}{2})\cap G_{n+1}$.

This construction implies that $\overline{B}_d(x_{n+1},r_{n+1})\subseteq\overline{B}_d(x_n,r_n)$ and $d(x_{n+1},x_n)<2^{-n}$. Therefore $(x_n)$ is Cauchy and converges to some $x\in X$. For each $n$, all terms $x_m$ with $m\geq n$ belong to $\overline{B}_d(x_n,r_n)$; since this ball is closed, $x\in\overline{B}_d(x_n,r_n)\subseteq G_n$. Moreover, $x\in\overline{B}_d(x_1,r_1)\subseteq U$. Consequently, $x\in U\cap\displaystyle\bigcap_{n=1}^{\infty}G_n$.
\end{proof}

\section{Normed spaces}\label{ap:espacios-normados}

Normed spaces provide the basic setting for measuring the size of vectors and operators. In this section we fix the notation for $\mathcal L(X,Y)$, recall the fundamental theorems on continuous linear operators, and introduce closed and compact operators. The exposition follows, among other references, \cite[Chapters 1 and 2]{DrabekMilota2013} and \cite{brezis_functional_2011,kesavan_functional_2023,Rudin1991}.

\begin{definition}\label{def:b2-espacios-normados-norma}\index{norm}
Let $X$ be a vector space over $\mathbb K=\mathbb R$ or $\mathbb C$. A \textbf{norm} on $X$ is a function $\|\cdot\|_X\colon X\longrightarrow[0,\infty)$ such that, for all $x,y\in X$ and $\lambda\in\mathbb K$,
\begin{enumerate}[label=(\alph*)]
\item $\|x\|_X=0$ if and only if $x=0$;
\item $\|\lambda x\|_X=|\lambda|\|x\|_X$;
\item $\|x+y\|_X\leq\|x\|_X+\|y\|_X$.
\end{enumerate}
The pair $(X,\|\cdot\|_X)$ is called a \textbf{normed space}. The distance induced by the norm is $d_X(x,y):=\|x-y\|_X$.
\end{definition}

\begin{definition}\label{def:b2-espacios-normados-espacio-de-banach}\index{Banach space@Banach space}
A normed space is called a \textbf{Banach space} if it is complete with respect to the distance induced by its norm.
\end{definition}

The metric completion of a normed space naturally retains the vector space operations.

\begin{theorem}[Completion of a normed space]\label{completación de un espacio normado}\index{completion of a normed space@completion of a normed space}
For every normed space $(X,\|\cdot\|_X)$, there exist a Banach space $(\widehat X,\|\cdot\|_{\widehat X})$ and a linear isometry $\iota\colon X\longrightarrow\widehat X$ with dense image. This completion is unique up to linear isometries that fix the copy of $X$.
\end{theorem}

\begin{proof}
Let $(\widehat X,\widehat d)$ be a completion of the metric space $(X,d_X)$, and let $\iota\colon X\longrightarrow\widehat X$ be the isometry with dense image provided by Theorem~\ref{completación de espacio métrico}. Given $u,v\in\widehat X$, choose sequences $(x_n)$ and $(y_n)$ in $X$ such that $\iota(x_n)\to u$ and $\iota(y_n)\to v$, and define
\[
u+v:=\lim_{n\to\infty}\iota(x_n+y_n),
\qquad
\lambda u:=\lim_{n\to\infty}\iota(\lambda x_n).
\]
The inequalities $\|(x_n+y_n)-(x_m+y_m)\|_X\leq\|x_n-x_m\|_X+\|y_n-y_m\|_X$ and $\|\lambda x_n-\lambda x_m\|_X=|\lambda|\|x_n-x_m\|_X$ show that the limits exist. The same estimates show that they do not depend on the chosen sequences. These operations make $\widehat X$ a vector space and make $\iota$ linear.

Define $\|u\|_{\widehat X}:=\widehat d(u,0)$. If $\iota(x_n)\to u$, then $\|u\|_{\widehat X}=\displaystyle\lim_{n\to\infty}\|x_n\|_X$. The norm properties follow by taking limits in the corresponding properties of $\|\cdot\|_X$. Since $\widehat X$ is complete with respect to $\widehat d$, it is a Banach space. Metric uniqueness of the completion and density of $\iota(X)$ imply that the isometry between two completions is necessarily linear.
\end{proof}

\subsection{Continuous linear operators}

\begin{definition}\label{def:b2-espacios-normados-operador-lineal}\index{linear operator}
Let $X$ and $Y$ be vector spaces over the same field $\mathbb K$. A map $T\colon X\longrightarrow Y$ is called a \textbf{linear operator} if $T(x+y)=T(x)+T(y)$ and $T(\lambda x)=\lambda T(x)$ for all $x,y\in X$ and $\lambda\in\mathbb K$.
\end{definition}

\begin{definition}\label{def:b2-espacios-normados-operador-acotado}\index{linear operator!bounded}
Let $X$ and $Y$ be normed spaces. A linear operator $T\colon X\longrightarrow Y$ is called \textbf{bounded} if there exists $C\geq0$ such that $\|T(x)\|_Y\leq C\|x\|_X$ for every $x\in X$.

Let $\mathcal L(X,Y)$ denote the vector space of all continuous linear operators from $X$ to $Y$, and write $\mathcal L(X):=\mathcal L(X,X)$. The \textbf{operator norm} is given by
\[
\|T\|_{\mathcal L(X,Y)}
:=\sup_{\substack{x\in X\\\|x\|_X\leq1}}\|T(x)\|_Y
=\sup_{x\in X\setminus\{0\}}\frac{\|T(x)\|_Y}{\|x\|_X}.
\]
When there is no risk of confusion, we will simply write $\|T\|$. The identity operator on $X$ is denoted by $I_X$, or by $I$ when the space is clear from context.
\end{definition}

\begin{definition}\label{def:dual-continuo-espacio-normado}\index{continuous dual}
The \textbf{continuous dual} of a normed space $X$ is $X':=\mathcal L(X,\mathbb K)$. Its elements are called \textbf{continuous linear functionals}.
\end{definition}

\begin{proposition}\label{prop:b2-espacios-normados-operador-lineal-dos-espacios-normados-funcion}
Let $X$ and $Y$ be normed spaces, and let $T\colon X\longrightarrow Y$ be linear. The following statements are equivalent:
\begin{enumerate}[label=(\alph*)]
\item $T$ is continuous on $X$;
\item $T$ is continuous at $0$;
\item $T$ is bounded.
\end{enumerate}
In this case, $\|T(x)\|_Y\leq\|T\|_{\mathcal L(X,Y)}\|x\|_X$ for every $x\in X$.
\end{proposition}

\begin{proof}
The implication $(a)\Rightarrow(b)$ is immediate. Suppose that $(b)$. There exists $\delta>0$ such that $\|T(x)\|_Y<1$ whenever $\|x\|_X<\delta$. If $x\neq0$, then $\left\|\displaystyle\frac{\delta x}{2\|x\|_X}\right\|_X=\displaystyle\frac{\delta}{2}$ and, by linearity, $\|T(x)\|_Y<\displaystyle\frac{2}{\delta}\|x\|_X$. The same inequality is trivial for $x=0$, so $T$ is bounded. Finally, if $(c)$ holds with a constant $C$, then $\|T(x)-T(y)\|_Y\leq C\|x-y\|_X$ for all $x,y\in X$, and hence $T$ is Lipschitz and continuous.
\end{proof}

\begin{proposition}\label{prop:propiedades-norma-operadores}
Let $X$, $Y$, and $Z$ be normed spaces. If $T,S\in\mathcal L(X,Y)$, $R\in\mathcal L(Y,Z)$, and $\lambda\in\mathbb K$, then
\[
\|T+S\|\leq\|T\|+\|S\|,
\qquad
\|\lambda T\|=|\lambda|\|T\|,
\qquad
\|R\circ T\|\leq\|R\|\|T\|.
\]
In particular, $\mathcal L(X)$ is a normed algebra under composition.
\end{proposition}

\begin{proof}
The first two identities follow directly from the definition. For the third, $\|R(Tx)\|_Z\leq\|R\|\|Tx\|_Y\leq\|R\|\|T\|\|x\|_X$ for every $x\in X$.
\end{proof}

\begin{theorem}\label{teo:completitud-espacio-operadores}
Let $X$ be a normed space and $Y$ a Banach space. Then $\mathcal L(X,Y)$ is a Banach space. In particular, $X'$ is a Banach space for every normed space $X$.
\end{theorem}

\begin{proof}
Let $(T_n)$ be a Cauchy sequence in $\mathcal L(X,Y)$. For each $x\in X$, the inequality $\|T_nx-T_mx\|_Y\leq\|T_n-T_m\|\|x\|_X$ shows that $(T_nx)$ is Cauchy in $Y$. Define $T(x):=\displaystyle\lim_{n\to\infty}T_nx$. The map $T$ is linear. Moreover, there exists $N\in\mathbb N$ such that $\|T_n-T_m\|<1$ for $n,m\geq N$, so $\|T_n\|\leq1+\|T_N\|$ for $n\geq N$. Taking limits gives $\|T(x)\|_Y\leq(1+\|T_N\|)\|x\|_X$, and hence $T\in\mathcal L(X,Y)$.

Given $\varepsilon>0$, there exists $N_\varepsilon\in\mathbb N$ such that $\|T_n-T_m\|<\varepsilon$ for $n,m\geq N_\varepsilon$. Letting $m\to\infty$, we obtain $\|(T_n-T)x\|_Y\leq\varepsilon\|x\|_X$ for every $x\in X$ and $n\geq N_\varepsilon$. Therefore, $\|T_n-T\|\leq\varepsilon$.
\end{proof}

\begin{theorem}[Neumann series]\label{teo:serie-de-neumann}\index{Neumann series}
Let $X$ be a Banach space and let $T\in\mathcal L(X)$ satisfy $\|T\|<1$. Then $I-T$ is invertible and
\[
(I-T)^{-1}=\displaystyle\sum_{n=0}^{\infty}T^n
\]
in the norm of $\mathcal L(X)$. Moreover, $\|(I-T)^{-1}\|\leq\displaystyle\frac{1}{1-\|T\|}$.
\end{theorem}

\begin{proof}
The series converges in $\mathcal L(X)$ because $\displaystyle\sum_{n=0}^{\infty}\|T^n\|\leq\displaystyle\sum_{n=0}^{\infty}\|T\|^n<\infty$ and $\mathcal L(X)$ is complete. If $S_N:=\displaystyle\sum_{n=0}^{N}T^n$, then $(I-T)S_N=S_N(I-T)=I-T^{N+1}$. Letting $N\to\infty$ gives $(I-T)S=S(I-T)=I$, where $S:=\displaystyle\sum_{n=0}^{\infty}T^n$. The norm estimate follows from the geometric series.
\end{proof}

\begin{corollary}\label{cor:estabilidad-invertibilidad-operadores}
Let $X$ be a Banach space, let $T\in\mathcal L(X)$ be invertible, and let $S\in\mathcal L(X)$ satisfy $\|T^{-1}(S-T)\|<1$. Then $S$ is invertible.
\end{corollary}

\begin{proof}
Since $S=T[I+T^{-1}(S-T)]$, the conclusion follows from Theorem~\ref{teo:serie-de-neumann} applied to $-T^{-1}(S-T)$.
\end{proof}

\begin{definition}\label{def:resolvente-espectro-operador-acotado}\index{spectrum!of a bounded operator}
Let $X$ be a complex Banach space and let $T\in\mathcal L(X)$. The \textbf{resolvent set} of $T$ is
\[
\rho(T):=\{\lambda\in\mathbb C\mid \lambda I-T\text{ is invertible in }\mathcal L(X)\},
\]
and the \textbf{spectrum} of $T$ is $\sigma(T):=\mathbb C\setminus\rho(T)$. For $\lambda\in\rho(T)$, the operator $R(\lambda,T):=(\lambda I-T)^{-1}$ is called the \textbf{resolvent} of $T$.
\end{definition}

\begin{proposition}\label{prop:propiedades-basicas-resolvente}
If $X$ is a complex Banach space and $T\in\mathcal L(X)$, then $\{\lambda\in\mathbb C\mid |\lambda|>\|T\|\}\subseteq\rho(T)$ and $\rho(T)$ is open. In particular, $\sigma(T)$ is compact and is contained in $\{\lambda\in\mathbb C\mid |\lambda|\leq\|T\|\}$.
\end{proposition}

\begin{proof}
If $|\lambda|>\|T\|$, then $\lambda I-T=\lambda\left(I-\displaystyle\frac{T}{\lambda}\right)$ and the Neumann series gives its inverse. If $\lambda_0\in\rho(T)$, we write
\[
\lambda I-T=(\lambda_0I-T)\bigl[I+(\lambda-\lambda_0)(\lambda_0I-T)^{-1}\bigr].
\]
The second factor is invertible when $|\lambda-\lambda_0|\|(\lambda_0I-T)^{-1}\|<1$, so $\rho(T)$ is open. The remaining assertions are immediate.
\end{proof}

\begin{theorem}\label{teo:espectro-no-vacio-operador-acotado}
If $X\neq\{0\}$ is a complex Banach space and $T\in\mathcal L(X)$, then $\sigma(T)$ is nonempty.
\end{theorem}

The proof uses Liouville's theorem applied to the scalar functions $\lambda\mapsto\varphi((\lambda I-T)^{-1}x)$, with $x\in X$ and $\varphi\in X'$, and can be found in \cite[Section 2.1]{DrabekMilota2013} or \cite{Rudin1991}.

\subsection{Fundamental principles of operator theory}

\begin{theorem}[Banach--Steinhaus]\label{teo: banach steinhaus}\index{Banach Steinhaus@Banach--Steinhaus}
Let $X$ be a Banach space, let $Y$ be a normed space, and let $\mathcal F\subseteq\mathcal L(X,Y)$. If $\displaystyle\sup_{T\in\mathcal F}\|T(x)\|_Y<\infty$ for every $x\in X$, then $\displaystyle\sup_{T\in\mathcal F}\|T\|_{\mathcal L(X,Y)}<\infty$.
\end{theorem}

\begin{proof}
For each $n\in\mathbb N$, define $E_n:=\{x\in X\mid \|T(x)\|_Y\leq n\text{ for every }T\in\mathcal F\}$. Each $E_n$ is closed and, by pointwise boundedness, $X=\displaystyle\bigcup_{n=1}^{\infty}E_n$. Theorem~\ref{teo:categoria-de-baire} gives $N\in\mathbb N$, $x_0\in X$, and $r>0$ such that $B_{d_X}(x_0,r)\subseteq E_N$, where $d_X$ is the metric induced by the norm of $X$.

If $\|h\|_X<r$, then $x_0+h,x_0\in E_N$ and $\|T(h)\|_Y\leq\|T(x_0+h)\|_Y+\|T(x_0)\|_Y\leq2N$ for every $T\in\mathcal F$. For $x\neq0$, we apply this estimate to $h:=\displaystyle\frac{r x}{2\|x\|_X}$ and obtain $\|T(x)\|_Y\leq\displaystyle\frac{4N}{r}\|x\|_X$. Therefore, $\|T\|\leq\displaystyle\frac{4N}{r}$ for every $T\in\mathcal F$.
\end{proof}

\begin{theorem}[Open mapping theorem]\label{teo:mapeo-abierto-banach}\index{open mapping theorem}
Let $X$ and $Y$ be Banach spaces and let $T\in\mathcal L(X,Y)$ be surjective. Then $T$ is an open map.
\end{theorem}

\begin{proof}
Put $\mathbb B_X:=B_{d_X}(0,1)$, where $d_X$ is the distance
induced by the norm. Surjectivity gives
$Y=\bigcup_{n\geq1}T(n\mathbb B_X)$, and hence also
$Y=\bigcup_{n\geq1}\overline{T(n\mathbb B_X)}$.
By Theorem~\ref{teo:categoria-de-baire}, there are $n_0\geq1$,
$y_0\in Y$ and $\delta>0$ such that
\[
 B_{d_Y}(y_0,\delta)\subseteq\overline{T(n_0\mathbb B_X)}.
\]
If $\|w\|_Y<\delta$, both $y_0+w$ and $y_0$ are limits of
images of vectors of norm less than $n_0$. Subtracting the two
sequences gives $w\in\overline{T(2n_0\mathbb B_X)}$.
Dividing by $2n_0$ yields
\[
 B_{d_Y}(0,r)\subseteq\overline{T(\mathbb B_X)},
 \qquad r:=\frac{\delta}{2n_0}.
\]
By homogeneity, for every $a>0$,
$B_{d_Y}(0,ar)\subseteq\overline{T(a\mathbb B_X)}$.

Let $\|y\|_Y<\displaystyle\frac{r}{2}$ and put $e_0:=y$. Inductively construct
$x_n\in X$ with $\|x_n\|_X<2^{-n}$ and
\[
 e_n:=y-T\left(\sum_{j=1}^n x_j\right),
 \qquad \|e_n\|_Y<r2^{-n-1}.
\]
Indeed, if $\|e_{n-1}\|_Y<r2^{-n}$, the preceding inclusion with
$a=2^{-n}$ allows us to approximate $e_{n-1}$ by $Tx_n$, with
$\|x_n\|_X<2^{-n}$, to within an error less than $r2^{-n-1}$.
Completeness of $X$ implies that $\sum_{n\geq1}x_n$ converges.
Its sum $x$ satisfies
\[
 \|x\|_X\leq\sum_{n\geq1}\|x_n\|_X
 \leq\|x_1\|_X+\sum_{n\geq2}2^{-n}<1.
\]
Continuity of $T$ and $e_n\to0$ give $Tx=y$. Thus
$B_{d_Y}(0,\displaystyle\frac{r}{2})\subseteq T(\mathbb B_X)$.

Finally, let $U\subseteq X$ be open and let $y=Tx\in T(U)$,
with $x\in U$. Choose $\varepsilon>0$ such that
$B_{d_X}(x,\varepsilon)\subseteq U$. Linearity and the inclusion
just proved give
\[
 y+B_{d_Y}(0,\varepsilon \frac{r}{2})
 \subseteq T\bigl(B_{d_X}(x,\varepsilon)\bigr)\subseteq T(U).
\]
Every point of $T(U)$ is interior; therefore $T(U)$ is open.
\end{proof}

\begin{corollary}[Bounded inverse theorem]\label{teo:inverso-acotado}\index{bounded inverse theorem}
Let $X$ and $Y$ be Banach spaces. If $T\in\mathcal L(X,Y)$ is bijective, then $T^{-1}\in\mathcal L(Y,X)$.
\end{corollary}

\begin{proof}
Theorem~\ref{teo:mapeo-abierto-banach} implies that the inverse operator $T^{-1}$ is continuous.
\end{proof}

\begin{proposition}\label{prop:criterio-rango-cerrado-inyectivo}
Let $X$ and $Y$ be Banach spaces and let $T\in\mathcal L(X,Y)$ be injective. Then $\operatorname{Ran}(T)$ is closed if and only if there exists $c>0$ such that $\|Tx\|_Y\geq c\|x\|_X$ for every $x\in X$.
\end{proposition}

\begin{proof}
If the estimate holds and $Tx_n\to y$ in $Y$, then $(x_n)$ is Cauchy and converges to some $x\in X$; by continuity, $Tx=y$, so the range is closed. Conversely, if $\operatorname{Ran}(T)$ is closed, it is a Banach space and $T\colon X\longrightarrow\operatorname{Ran}(T)$ is bijective. Corollary~\ref{teo:inverso-acotado} implies that $T^{-1}$ is continuous, and therefore $\|x\|_X\leq\|T^{-1}\|\|Tx\|_Y$.
\end{proof}

Quotient spaces allow us to factor out a closed subspace while retaining a
normed structure and, when the original space is Banach, completeness.
Closedness is essential: without it, the expression below defines only
a seminorm.

\begin{definition}
\label{def:espacio-cociente-vectorial}
\index{quotient space!vector}
Let $X$ be a vector space and let $Y\subseteq X$ be a subspace. The
\textbf{quotient vector space} $X/Y$ is the set of classes
\[
X/Y:=\{x+Y\mid x\in X\},
\qquad
[x]:=x+Y,
\]
with the operations
\[
[x]+[z]:=[x+z],
\qquad
\alpha[x]:=[\alpha x]
\quad (\alpha\in\mathbb K).
\]
These operations do not depend on the representatives: if $x-x'\in Y$ and
$z-z'\in Y$, then $(x+z)-(x'+z')\in Y$ and
$\alpha x-\alpha x'\in Y$.
\end{definition}

\begin{proposition}[Quotient norm]
\label{prop:norma-cociente-banach}
\index{norm!quotient}
\index{quotient space!of a Banach space}
Let $X$ be a normed space and let $Y\subseteq X$ be a closed subspace. For
$[x]=x+Y\in X/Y$, set
\begin{equation}
\label{eq:norma-cociente-banach}
\|[x]\|_{X/Y}
:=
\inf_{y\in Y}\|x-y\|_X
=
\operatorname{dist}_X(x,Y).
\end{equation}
Then \eqref{eq:norma-cociente-banach} defines a norm on $X/Y$. The
canonical projection
\[
\pi\colon X\longrightarrow X/Y,
\qquad
\pi(x)=[x],
\]
is linear, continuous, and open, and satisfies $\|\pi\|\leq1$. If $X$ is a
Banach space, then $X/Y$ is also a Banach space.
\end{proposition}

\begin{proof}
If $x'=x+y_0$ with $y_0\in Y$, then
\[
\inf_{y\in Y}\|x'-y\|_X
=
\inf_{y\in Y}\|x-(y-y_0)\|_X
=
\inf_{z\in Y}\|x-z\|_X,
\]
so the expression does not depend on the representative of the class.
Homogeneity is immediate. Given $x,z\in X$ and $y_1,y_2\in Y$, we have
\[
\operatorname{dist}_X(x+z,Y)
\leq
\|(x+z)-(y_1+y_2)\|_X
\leq
\|x-y_1\|_X+\|z-y_2\|_X.
\]
Taking the infimum over $y_1$ and $y_2$ gives the triangle
inequality. Moreover,
$\|[x]\|_{X/Y}=0$ if and only if $x\in\overline Y=Y$, that is, if and only if
$[x]=0$. Therefore, \eqref{eq:norma-cociente-banach} is a norm.

The linearity of $\pi$ follows from the quotient operations, and
$\|\pi(x)\|_{X/Y}\leq\|x\|_X$, so $\pi$ is continuous and
$\|\pi\|\leq1$. If $r>0$, then for every class $[z]$ with
$\|[z]\|_{X/Y}<r$ there exists $y\in Y$ such that $\|z-y\|_X<r$ and
$[z]=[z-y]$. Consequently, every class of norm less than $r$ belongs to
$\pi(B_{d_X}(0,r))$. The reverse inclusion follows from
$\|\pi(x)\|_{X/Y}\leq\|x\|_X$. Therefore, the image of the open ball of
radius $r$ in $X$ is the open ball of radius $r$ in $X/Y$. Translations and
linearity show that $\pi$ is open.

Finally, suppose that $X$ is complete and let $([x_n])$ be a Cauchy sequence
in $X/Y$. Choose a subsequence $([x_{n_j}])$ such that
\[
\|[x_{n_{j+1}}]-[x_{n_j}]\|_{X/Y}<2^{-j-1}
\qquad (j\in\mathbb N).
\]
By the definition of the infimum, for each $j$ there exists $y_j\in Y$ such that
\[
\|(x_{n_{j+1}}-x_{n_j})-y_j\|_X<2^{-j}.
\]
Define $z_1:=x_{n_1}$ and, for $j\geq2$,
\[
z_j:=x_{n_j}-\sum_{i=1}^{j-1}y_i.
\]
Then $[z_j]=[x_{n_j}]$ and
\[
\|z_{j+1}-z_j\|_X
=
\|(x_{n_{j+1}}-x_{n_j})-y_j\|_X
<2^{-j}.
\]
Consequently, $(z_j)$ is Cauchy in $X$ and converges to some $z\in X$.
The continuity of $\pi$ implies that $[x_{n_j}]=[z_j]\to[z]$. Since the original sequence
is Cauchy, the entire sequence converges to $[z]$. This proves the completeness of
$X/Y$.
\end{proof}

Differential operators are usually unbounded when considered on a single function space. To handle them, we must distinguish the domain of the operator from the ambient space.

\begin{definition}\label{def:operador-densamente-definido-cerrado}\index{closed operator}\index{densely defined operator}
Let $X$ and $Y$ be normed spaces. A \textbf{linear operator} from $X$ to $Y$ is a linear map $A\colon \mathcal D(A)\subseteq X\longrightarrow Y$, where $\mathcal D(A)$ is a vector subspace called the \textbf{domain} of $A$; we also write $\operatorname{Dom}(A)$. The operator is \textbf{densely defined} if $\overline{\mathcal D(A)}=X$.

The \textbf{graph} of $A$ is $\mathcal G(A):=\{(x,A x)\mid x\in\mathcal D(A)\}\subseteq X\times Y$. We say that $A$ is \textbf{closed} if $\mathcal G(A)$ is closed. Equivalently, if $x_n\in\mathcal D(A)$, $x_n\to x$ in $X$, and $Ax_n\to y$ in $Y$, then $x\in\mathcal D(A)$ and $Ax=y$.
\end{definition}

\begin{theorem}[Closed graph theorem]\label{teo:grafica-cerrada}\index{closed graph theorem}
Let $X$ and $Y$ be Banach spaces and let $A\colon X\longrightarrow Y$ be linear. Then $A$ is continuous if and only if its graph is closed.
\end{theorem}

\begin{proof}
If $A$ is continuous, its graph is immediately seen to be closed. Conversely, suppose that $\mathcal G(A)$ is closed. Then $\mathcal G(A)$ is a Banach space with the norm inherited from $X\times Y$. The projection $\pi_1\colon \mathcal G(A)\longrightarrow X$, given by $\pi_1(x,Ax)=x$, is linear, continuous, and bijective. By Corollary~\ref{teo:inverso-acotado}, $\pi_1^{-1}$ is continuous. Since $A=\pi_2\circ\pi_1^{-1}$, where $\pi_2\colon X\times Y\longrightarrow Y$ is the second projection, $A$ is continuous.
\end{proof}

\begin{proposition}\label{prop:norma-grafica-operador-cerrado}
Let $X$ and $Y$ be Banach spaces and let $A\colon \mathcal D(A)\subseteq X\longrightarrow Y$ be linear. The formula $\|x\|_A:=\|x\|_X+\|Ax\|_Y$ defines a norm on $\mathcal D(A)$. The operator $A$ is closed if and only if $(\mathcal D(A),\|\cdot\|_A)$ is a Banach space.
\end{proposition}

\begin{proof}
The map $J\colon \mathcal D(A)\longrightarrow X\times Y$, given by $J(x)=(x,Ax)$, is an isometry from $(\mathcal D(A),\|\cdot\|_A)$ onto $\mathcal G(A)$ when $X\times Y$ is equipped with the norm $\|(x,y)\|:=\|x\|_X+\|y\|_Y$. Therefore, $\mathcal D(A)$ is complete if and only if $\mathcal G(A)$ is closed in the Banach space $X\times Y$.
\end{proof}

\begin{definition}\label{def:operador-cerrable-cerradura}\index{closable operator}
A linear operator $A\colon \mathcal D(A)\subseteq X\longrightarrow Y$ is called \textbf{closable} if the closure of $\mathcal G(A)$ is the graph of an operator. This operator is denoted by $\overline A$ and is called the \textbf{closure} of $A$.
\end{definition}

\begin{proposition}\label{prop:criterio-operador-cerrable}
A linear operator $A\colon \mathcal D(A)\subseteq X\longrightarrow Y$ is closable if and only if, for every sequence $(x_n)$ in $\mathcal D(A)$, the convergences $x_n\to0$ in $X$ and $Ax_n\to y$ in $Y$ imply $y=0$.
\end{proposition}

\begin{proof}
The closure of $\mathcal G(A)$ is the graph of an operator if and only if it contains no element $(0,y)$ with $y\neq0$. Since $X\times Y$ is metrizable, a point belongs to the closure if and only if it is the limit of a sequence of elements of $\mathcal G(A)$, which gives precisely the stated condition.
\end{proof}

\subsection{The Hahn--Banach theorem and duality}

The Hahn--Banach theorem allows us to extend linear functionals while retaining the estimates that ensure their continuity. Its consequences will be used both in the dual characterization of the norm and in separation arguments.

\begin{remark}\label{rem: sublineal}
A function $p\colon X\longrightarrow\mathbb R$ defined on a real vector space is called \textbf{sublinear} if $p(x+y)\leq p(x)+p(y)$ for all $x,y\in X$ and $p(\lambda x)=\lambda p(x)$ for every $x\in X$ and $\lambda\geq0$.
\end{remark}

\begin{theorem}[Hahn--Banach theorem, analytic form]\label{teo: hahn--banach, forma analitica}\index{Hahn Banach theorem@Hahn--Banach theorem!analytic form}
Let $X$ be a real vector space, let $M\subseteq X$ be a subspace, and let $p\colon X\longrightarrow\mathbb R$ be sublinear. If $\ell\colon M\longrightarrow\mathbb R$ is linear and $\ell(x)\leq p(x)$ for every $x\in M$, then there exists a linear functional $\widetilde\ell\colon X\longrightarrow\mathbb R$ such that $\widetilde\ell\restriction_M=\ell$ and $\widetilde\ell(x)\leq p(x)$ for every $x\in X$.
\end{theorem}
\begin{proof}
We begin with extension by one dimension. Let $N\subsetneq X$ be a subspace, let $f\colon N\to\mathbb R$ be linear and dominated by $p$, and let $x_0\in X\setminus N$. Define
\[
\alpha:=
\sup_{y\in N}\bigl(f(y)-p(y-x_0)\bigr),
\qquad
\beta:=
\inf_{z\in N}\bigl(p(z+x_0)-f(z)\bigr).
\]
For any $y,z\in N$, the linearity of $f$, domination by $p$, and subadditivity of $p$ imply
\[
f(y)+f(z)
=f(y+z)
\leq p(y+z)
\leq p(y-x_0)+p(z+x_0).
\]
Therefore,
\[
f(y)-p(y-x_0)\leq p(z+x_0)-f(z)
\]
for every pair $y,z\in N$. In particular, $\alpha\leq\beta$, and we may choose $a\in[\alpha,\beta]$.

Every element of $N+\mathbb R x_0$ has a unique representation as $y+t x_0$, with $y\in N$ and $t\in\mathbb R$. Define
\[
f_1(y+t x_0):=f(y)+ta.
\]
The map $f_1$ is linear and extends $f$. Let us verify domination. If $t>0$, then $a\leq\beta$ and, taking $z=\displaystyle\frac{y}{t}$, we obtain
\[
f_1(y+t x_0)
=t\bigl(f(\frac{y}{t})+a\bigr)
\leq t\,p(\frac{y}{t}+x_0)
=p(y+t x_0).
\]
If $t<0$, write $t=-s$ with $s>0$. Since $a\geq\alpha$, applied now to $\displaystyle\frac{y}{s}$, we have
\[
f_1(y-sx_0)
=s\bigl(f(\frac{y}{s})-a\bigr)
\leq s\,p(\frac{y}{s}-x_0)
=p(y-sx_0).
\]
The case $t=0$ is precisely the hypothesis $f\leq p$ on $N$. Thus, $f_1$ is dominated by $p$ on all of $N+\mathbb R x_0$.

Consider the set \(\mathcal E\) of all pairs $(N,f)$ such that $M\subseteq N\subseteq X$, $f\colon N\to\mathbb R$ is linear, $f\restriction_M=\ell$, and $f\leq p$ on $N$. Order this set by extension:
\[
(N_1,f_1)\preceq(N_2,f_2)
\quad\Longleftrightarrow\quad
N_1\subseteq N_2\ \text{and}\ f_2\restriction_{N_1}=f_1.
\]
The set is nonempty because it contains $(M,\ell)$. If \(\mathcal C\subseteq\mathcal E\) is a chain, the union of its domains
\[
N_{\mathcal C}:=\bigcup_{(N,f)\in\mathcal C}N
\]
is a subspace. On this subspace, define \(f_{\mathcal C}(x):=f(x)\) using any element of the chain whose domain contains $x$. The definition is independent of this choice, since any two elements of a chain are comparable; moreover, \(f_{\mathcal C}\) is linear, extends \(\ell\), and is dominated by \(p\). Thus every chain has an upper bound in \(\mathcal E\).

Zorn's lemma, Theorem~\ref{teo:lema-zorn}, gives a maximal element $(N_{\displaystyle\max},f_{\displaystyle\max})$. If $N_{\displaystyle\max}\neq X$, we may choose $x_0\in X\setminus N_{\displaystyle\max}$ and apply the preceding one-dimensional construction, obtaining a proper extension of $f_{\displaystyle\max}$ that is still dominated by $p$. This contradicts maximality. Hence $N_{\displaystyle\max}=X$, and \(\widetilde\ell:=f_{\displaystyle\max}\) is the desired extension.
\end{proof}

\begin{remark}\label{rem:hahn-banach-complejo}
On a complex vector space, the corresponding version is obtained by assuming that $p$ is subadditive, that $p(\lambda x)=|\lambda|p(x)$ for all $\lambda\in\mathbb C$ and $x\in X$, and that $|\ell(x)|\leq p(x)$ on the initial subspace. Indeed, regard \(X\) and \(M\) as real vector spaces and apply the theorem to \(u:=\operatorname{Re}\ell\). We obtain a real-linear extension \(U\colon X\to\mathbb R\) such that \(U\leq p\). Since \(p(-x)=p(x)\), we also have
\[
-U(x)=U(-x)\leq p(-x)=p(x),
\]
so \(|U(x)|\leq p(x)\). Define
\[
\widetilde\ell(x):=U(x)-iU(ix).
\]
The real linearity of \(U\) shows that
\[
\widetilde\ell(ix)
=U(ix)-iU(-x)
=i\widetilde\ell(x),
\]
and therefore \(\widetilde\ell\) is complex-linear. If \(x\in M\), then
\[
U(x)=\operatorname{Re}\ell(x),
\qquad
U(ix)=\operatorname{Re}(i\ell(x))=-\operatorname{Im}\ell(x),
\]
so \(\widetilde\ell(x)=\ell(x)\).

Finally, given \(x\in X\), choose \(\lambda\in\mathbb C\), \(|\lambda|=1\), so that \(\lambda\widetilde\ell(x)=|\widetilde\ell(x)|\). Since \(\operatorname{Re}\widetilde\ell=U\),
\[
|\widetilde\ell(x)|
=\operatorname{Re}\widetilde\ell(\lambda x)
=U(\lambda x)
\leq p(\lambda x)
=p(x).
\]
This proves the complex version and, in particular, the estimate used in the following corollary.
\end{remark}

\begin{corollary}[Norm-preserving extension]\label{cor:extension-hahn-banach-preserva-norma}\index{Hahn Banach theorem@Hahn--Banach theorem!norm-preserving extension}
Let $X$ be a real or complex normed space, let $M\subseteq X$ be a subspace equipped with the induced norm, and let $\ell\in M'$. Then there exists $\widetilde\ell\in X'$ such that $\widetilde\ell\restriction_M=\ell$ and $\|\widetilde\ell\|_{X'}=\|\ell\|_{M'}$.
\end{corollary}

\begin{proof}
Define $p(x):=\|\ell\|_{M'}\|x\|_X$. In the real case, $\ell(x)\leq|\ell(x)|\leq p(x)$ for every $x\in M$, so Theorem~\ref{teo: hahn--banach, forma analitica} gives an extension $\widetilde\ell$ dominated by $p$. In the complex case, we apply Remark~\ref{rem:hahn-banach-complejo}. In both cases, $|\widetilde\ell(x)|\leq\|\ell\|_{M'}\|x\|_X$, so $\|\widetilde\ell\|_{X'}\leq\|\ell\|_{M'}$. The reverse inequality follows from $\widetilde\ell\restriction_M=\ell$.
\end{proof}

\begin{corollary}[Dual characterization of the norm]\label{cor:caracterizacion-dual-norma}\index{norm!dual characterization}\index{continuous dual!characterization of the norm}
Let $X$ be a real or complex normed space. For each $x\in X$,
\[
\|x\|_X
=
\max\left\{|\varphi(x)|\middle|\varphi\in X',\ \|\varphi\|_{X'}\leq1\right\}.
\]
More precisely, if $x\neq0$, there exists $\varphi\in X'$ such that $\|\varphi\|_{X'}=1$ and $\varphi(x)=\|x\|_X$.
\end{corollary}

\begin{proof}
If $x\neq0$, define $\ell\colon \operatorname{span}\{x\}\longrightarrow\mathbb K$ by $\ell(\lambda x):=\lambda\|x\|_X$. Then $\|\ell\|=1$. Corollary~\ref{cor:extension-hahn-banach-preserva-norma} gives $\varphi\in X'$ with $\|\varphi\|=1$ and $\varphi(x)=\|x\|_X$. On the other hand, $|\psi(x)|\leq\|\psi\|\|x\|_X\leq\|x\|_X$ for every $\psi\in X'$ with $\|\psi\|\leq1$. The case $x=0$ is immediate.
\end{proof}

\subsection{Adjoint and compact operators on Banach spaces}

For compactness of the adjoint and the Riesz--Schauder alternative, we follow
\cite[Proposition~2.2.8 and Theorem~2.2.9]{DrabekMilota2013}.

\begin{definition}\label{def:adjunto-banach-operador-acotado}\index{adjoint operator!on Banach spaces}
Let $X$ and $Y$ be normed spaces and let $T\in\mathcal L(X,Y)$. The \textbf{Banach adjoint} of $T$ is the operator $T'\colon Y'\longrightarrow X'$ defined by $T'(\varphi):=\varphi\circ T$.
\end{definition}

\begin{proposition}\label{prop:propiedades-adjunto-banach}
Let $X$, $Y$, and $Z$ be normed spaces, and let $T\in\mathcal L(X,Y)$ and $S\in\mathcal L(Y,Z)$. Then $T'\in\mathcal L(Y',X')$, $\|T'\|=\|T\|$, and $(S\circ T)'=T'\circ S'$.
\end{proposition}

\begin{proof}
Composition with $T$ is linear in the functional. For $\varphi\in Y'$
and $x\in X$,
\[
 |(T'\varphi)(x)|=|\varphi(Tx)|
 \leq\|\varphi\|_{Y'}\|T\|\|x\|_X,
\]
so $T'\in\mathcal L(Y',X')$ and $\|T'\|\leq\|T\|$.
By Corollary~\ref{cor:caracterizacion-dual-norma}, for every
$x\in X$ with $\|x\|_X\leq1$ there is $\varphi\in Y'$ with
$\|\varphi\|_{Y'}\leq1$ and $|\varphi(Tx)|=\|Tx\|_Y$.
The choice $\varphi=0$ applies when $Tx=0$, including $Y=\{0\}$.
Then
\[
 \|Tx\|_Y=|(T'\varphi)(x)|\leq\|T'\|.
\]
Taking the supremum over the unit ball of $X$ gives $\|T\|\leq\|T'\|$.
Finally, for $\psi\in Z'$ and $x\in X$,
\[
 ((S\circ T)'\psi)(x)=\psi(S(Tx))
 =(T'(S'\psi))(x),
\]
which proves $(S\circ T)'=T'\circ S'$.
\end{proof}

\begin{definition}\label{def:b2-espacios-normados-operador-compacto}\index{compact operator}
Let $X$ and $Y$ be normed spaces. A linear operator $K\colon X\longrightarrow Y$ is called \textbf{compact} if, for every bounded sequence $(x_n)$ in $X$, the sequence $(Kx_n)$ has a subsequence converging in $Y$. We denote the space of compact operators from $X$ to $Y$ by $\mathcal K(X,Y)$ and write $\mathcal K(X):=\mathcal K(X,X)$.
\end{definition}

\begin{proposition}\label{prop:b2-espacios-normados-operador-lineal-compacto-operador-acotado-decir}
Every compact linear operator is bounded. Moreover, every \emph{continuous}
linear operator of finite rank is compact.
\end{proposition}

\begin{proof}
If a compact operator $K$ were unbounded, there would exist a sequence
$(x_n)$ with $\|x_n\|_X\leq1$ and $\|Kx_n\|_Y\geq n$, which would prevent
$(Kx_n)$ from having a convergent subsequence. If $K$ is continuous and has
finite rank, the image of every bounded sequence is bounded in the
finite-dimensional space $\operatorname{Ran}(K)$ and therefore has a
convergent subsequence.
\end{proof}

\begin{proposition}\label{prop:propiedades-operadores-compactos}
Let $X$, $Y$, and $Z$ be normed spaces.
\begin{enumerate}[label=(\alph*)]
\item If $K\in\mathcal K(X,Y)$ and $S\in\mathcal L(Y,Z)$, then $S\circ K\in\mathcal K(X,Z)$.
\item If $T\in\mathcal L(Z,X)$ and $K\in\mathcal K(X,Y)$, then $K\circ T\in\mathcal K(Z,Y)$.
\item If $Y$ is Banach, then $\mathcal K(X,Y)$ is a closed subspace of $\mathcal L(X,Y)$.
\end{enumerate}
\end{proposition}

\begin{proof}
Let $(x_j)$ be bounded in $X$. If $K$ is compact, a subsequence
of $(Kx_j)$ converges in $Y$, and continuity of $S$ gives convergence
of its images under $S$. This proves (a).
For (b), a bounded sequence $(z_j)$ in $Z$ has bounded images
$(Tz_j)$ in $X$, so compactness of $K$ gives a convergent
subsequence of $(KTz_j)$.

To prove that $\mathcal K(X,Y)$ is a subspace, let $K,L$ be compact.
First extract a subsequence along which $Kx_j$ converges, and then
one of it along which $Lx_j$ converges. Along the latter,
$(K+L)x_j$ also converges. Scalar multiplication preserves
convergence, and the zero operator is compact. Thus linear
combinations are compact.

Now suppose $Y$ is Banach and $K_n\to K$ in operator norm,
with each $K_n$ compact. Given a bounded sequence $(x_j)$,
successively extract nested subsequences along which
$K_1x_j,K_2x_j,\ldots$, respectively, converge.
Choose increasing indices $j_k$, with $j_k$ belonging to the
$k$-th subsequence. For fixed $n$, $(x_{j_k})_{k\geq n}$ is a
subsequence of the $n$-th extraction, so $(K_nx_{j_k})_k$
converges. Put $M:=\sup_j\|x_j\|_X<\infty$. Then
\[
 \|Kx_{j_k}-Kx_{j_l}\|_Y
 \leq2M\|K-K_n\|+\|K_nx_{j_k}-K_nx_{j_l}\|_Y.
\]
If $M=0$, all images are zero. Otherwise, given $\varepsilon>0$,
choose $n$ so that the first term is less than $\displaystyle\frac{\varepsilon}{2}$,
and then choose $k,l$ large enough for the second to be less than
$\displaystyle\frac{\varepsilon}{2}$. The images under $K$ are Cauchy and converge
by completeness of $Y$. Hence $K$ is compact and
$\mathcal K(X,Y)$ is closed.
\end{proof}

\begin{theorem}[Continuous and compact extension from a dense subspace]
\label{teo:extension-operadores-subespacio-denso}
\index{extension of linear operators!from a dense subspace}
Let $X$ be a normed space, let $D\subseteq X$ be a dense vector subspace, and let $Y$ be a Banach space. Every continuous linear operator $T\colon D\longrightarrow Y$ has a unique continuous linear extension $\widetilde T\colon X\longrightarrow Y$. Moreover,
\[
\|\widetilde T\|_{\mathcal L(X,Y)}=\|T\|_{\mathcal L(D,Y)}.
\]
If $T$ is compact, then $\widetilde T$ is also compact.
\end{theorem}

\begin{proof}
Let $x\in X$ and choose a sequence $(x_n)\subseteq D$ such that $x_n\to x$. The estimate $\|T(x_n)-T(x_m)\|_Y\leq\|T\|\|x_n-x_m\|_X$ shows that $(T(x_n))$ is Cauchy in $Y$. Since $Y$ is complete, we may define $\widetilde T(x):=\displaystyle\lim_{n\to\infty}T(x_n)$. If $(y_n)\subseteq D$ also converges to $x$, then $\|T(x_n)-T(y_n)\|_Y\leq\|T\|\|x_n-y_n\|_X\to0$, so the limit does not depend on the chosen sequence.

Linearity follows by approximating two points of $X$ by sequences in $D$ and passing to the limit. Similarly,
\[
\|\widetilde T(x)\|_Y
=\lim_{n\to\infty}\|T(x_n)\|_Y
\leq\|T\|\lim_{n\to\infty}\|x_n\|_X
=\|T\|\|x\|_X,
\]
so $\widetilde T$ is continuous and $\|\widetilde T\|\leq\|T\|$. The reverse inequality follows from $\widetilde T\restriction_D=T$. Uniqueness follows from the density of $D$: two continuous extensions that agree on $D$ agree on all of $X$.

Now suppose that $T$ is compact and let $(x_n)$ be a bounded sequence in $X$. For each $n$, choose $y_n\in D$ such that $\|x_n-y_n\|_X<\displaystyle\frac{1}{n}$. The sequence $(y_n)$ is also bounded, so there exists a subsequence $(y_{n_k})$ for which $(T(y_{n_k}))$ converges in $Y$. On the other hand,
\[
\|\widetilde T(x_{n_k})-T(y_{n_k})\|_Y
\leq
\|T\|\|x_{n_k}-y_{n_k}\|_X
<\frac{\|T\|}{n_k},
\]
and the right-hand side tends to zero. Consequently, $(\widetilde T(x_{n_k}))$ converges in $Y$, proving that $\widetilde T$ is compact.
\end{proof}

\begin{corollary}[Extension to the completions]
\label{cor:extension-operadores-completaciones}
\index{completion!extension of operators}
Let $X$ and $Y$ be normed spaces, and let $\widehat X$ and $\widehat Y$ be their completions. Every operator $T\in\mathcal L(X,Y)$ has a unique extension $\widehat T\in\mathcal L(\widehat X,\widehat Y)$ compatible with the canonical inclusions. We have $\|\widehat T\|=\|T\|$ and, if $T$ is compact, then $\widehat T$ is compact.
\end{corollary}

\begin{proof}
Identify $X$ and $Y$ with their isometric images in $\widehat X$ and $\widehat Y$. The result follows by applying Theorem~\ref{teo:extension-operadores-subespacio-denso} to the operator $X\longrightarrow\widehat Y$ induced by $T$.
\end{proof}

\begin{theorem}[Schauder's theorem]\label{teo:schauder-adjunto-compacto}\index{Schauder's theorem!adjoint of a compact operator}
Let $X$ and $Y$ be Banach spaces and let $K\in\mathcal L(X,Y)$. Then $K$ is compact if and only if $K'\colon Y'\longrightarrow X'$ is compact.
\end{theorem}

\begin{proof}
First suppose that $K$ is compact and set
\[
 C:=\overline{\{Kx\mid \|x\|_X\leq1\}}\subseteq Y.
\]
Every sequence in $C$ has a convergent subsequence. Indeed, approximate
its term $j$ by an element $Kx_j$, with $\|x_j\|_X\leq1$, at a distance
less than $\displaystyle\frac{1}{j}$. The compactness of $K$ gives a convergent subsequence
of $(Kx_j)$ and hence of the original sequence, whose limit
belongs to $C$ because this set is closed. In particular, $C$ has a finite net
of any radius $\delta>0$: if no finite union of balls of that
radius covered $C$, we could inductively choose a sequence of points
whose pairwise distances are at least $\delta$, contradicting the preceding
property.

Let $(\varphi_j)$ be a bounded sequence in $Y'$ and choose $M\geq1$ such that
$\|\varphi_j\|_{Y'}\leq M$ for every $j$. The restrictions
$f_j:=\varphi_j\restriction_C$ satisfy
\[
 |f_j(y)|\leq M\|K\|,
 \qquad
 |f_j(y)-f_j(z)|\leq M\|y-z\|_Y,
 \qquad y,z\in C.
\]
Let us give the uniform compactness argument needed here. For each
integer $q\geq1$, choose a finite net of radius $\displaystyle\frac{1}{q}$ in $C$. The union
of these nets is a countable dense subset of $C$. A diagonal argument,
applied to the bounded scalar sequences at its points, gives a
subsequence of $(f_j)$ that converges at every point of that union. Retain
the same indices for the subsequence.

Given $\varepsilon>0$, choose from these points a finite net
$\{y_1,\ldots,y_l\}$ of radius less than $\displaystyle\frac{\varepsilon}{4M}$. For sufficiently large
$j,k$,
$|f_j(y_i)-f_k(y_i)|<\displaystyle\frac{\varepsilon}{2}$ at every point of the net.
If $y\in C$ and $\|y-y_i\|_Y<\displaystyle\frac{\varepsilon}{4M}$, then
\[
 |f_j(y)-f_k(y)|
 \leq2M\|y-y_i\|_Y+|f_j(y_i)-f_k(y_i)|<\varepsilon.
\]
Thus the subsequence is uniformly Cauchy on $C$. Consequently,
\[
 \|K'\varphi_j-K'\varphi_k\|_{X'}
 =\sup_{\substack{x\in X\\\|x\|_X\leq1}}|(\varphi_j-\varphi_k)(Kx)|
 \leq\sup_{y\in C}|f_j(y)-f_k(y)|\longrightarrow0.
\]
Theorem~\ref{teo:completitud-espacio-operadores} shows that $X'$ is
complete; the subsequence of $(K'\varphi_j)$ converges in $X'$. This proves
that $K'$ is compact.

Conversely, suppose that $K'$ is compact. The duals $Y'$ and $X'$
are Banach spaces, so the implication already proved, applied to
$K'$, shows that $K''\colon X''\to Y''$ is compact. For a normed space
$Z$, define
\[
 J_Z\colon Z\longrightarrow Z'',
 \qquad (J_Zz)(\ell):=\ell(z),\quad \ell\in Z'.
\]
Corollary~\ref{cor:caracterizacion-dual-norma} gives
$\|J_Zz\|_{Z''}=\|z\|_Z$, so $J_Z$ is a linear isometry.
Moreover, evaluating at $\varphi\in Y'$ gives
\[
 (K''J_Xx)(\varphi)
 =(J_Xx)(K'\varphi)
 =\varphi(Kx)
 =(J_YKx)(\varphi).
\]
Hence $K''J_X=J_YK$. If $(x_j)$ is bounded in $X$, the compactness of
$K''$ gives a subsequence along which $(J_YKx_j)$ converges in
$Y''$. Since $J_Y$ is an isometry, the corresponding subsequence of
$(Kx_j)$ is Cauchy in $Y$ and converges because $Y$ is complete. Thus $K$
is compact.
\end{proof}

\begin{theorem}[Riesz--Schauder theorem]\label{teo:riesz-schauder}\index{Riesz--Schauder theorem}
Let $X$ be a Banach space and let $K\in\mathcal K(X)$. Then $\ker(I-K)$ is finite-dimensional, $\operatorname{Ran}(I-K)$ is closed, and
\[
\operatorname{Ran}(I-K)=X
\quad\Longleftrightarrow\quad
\ker(I-K)=\{0\}.
\]
Moreover, $\dim\ker(I-K)=\dim\ker(I-K')$.
\end{theorem}

\begin{proof}
Set $T:=I-K$ and $N:=\ker T$. We first use an elementary separation
fact. If $F$ is a proper closed subspace of a normed space
$E$, choose $z\in E\setminus F$ and set $d:=\operatorname{dist}(z,F)>0$.
There exists $y\in F$ such that $\|z-y\|<2d$. The vector
\[
 u:=\frac{z-y}{\|z-y\|}
\]
satisfies $\|u\|=1$ and
$\operatorname{dist}(u,F)=\displaystyle\frac{d}{\|z-y\|}>\displaystyle\frac{1}{2}$.
Recall also that every finite-dimensional subspace is closed and
that its coordinate functionals are continuous: for a basis
$e_1,\ldots,e_l$, the function
$\displaystyle \|\displaystyle\sum_{i=1}^{l} a_i e_i\|$ attains a positive minimum on the Euclidean sphere
$\displaystyle \sum_{i=1}^{l}|a_i|^2=1$. Homogeneity gives equivalence with the coordinate
norm, from which both assertions follow.

The subspace $N$ is closed because $T$ is continuous. If it were
infinite-dimensional, we could apply the separation fact successively to the
subspaces spanned by the vectors already chosen. We would obtain a
sequence $(u_j)\subseteq N$, with $\|u_j\|_X=1$, such that
$\|u_j-u_k\|_X>\displaystyle\frac{1}{2}$ for $j\ne k$. But $Ku_j=u_j$, contradicting
the compactness of $K$. Thus, $\dim N<\infty$.

We now prove the estimate
\[
 \operatorname{dist}_X(x,N)\leq C\|Tx\|_X,
 \qquad x\in X,
\]
for some constant $C>0$. If it were false, after normalizing and
subtracting an element of $N$, we could choose $x_j\in X$ such that
\[
 \operatorname{dist}_X(x_j,N)=1,
 \qquad \|x_j\|_X<2,
 \qquad \|Tx_j\|_X<\frac{1}{j}.
\]
In detail, take a vector whose distance to $N$ is greater than
$j\|Tx\|_X$, divide by this distance, and choose a representative of
the same class modulo $N$ with norm less than two. The compactness of $K$
allows us to pass to a subsequence along which $Kx_j$ converges. Since
$x_j=Tx_j+Kx_j$, $x_j$ also converges, say to $x$. By continuity,
$Tx=0$ and $\operatorname{dist}_X(x,N)=1$, a contradiction.
By Proposition~\ref{prop:norma-cociente-banach},
the estimate is a lower bound for the
operator induced by $T$ on $X/N$.

Now let $Tx_j\to y$ in $X$. Choose $n_j\in N$ so that
\[
 z_j:=x_j-n_j,
 \qquad
 \|z_j\|_X\leq2\operatorname{dist}_X(x_j,N)+\frac{1}{j}
 \leq2C\|Tx_j\|_X+\frac{1}{j}.
\]
The sequence $(z_j)$ is bounded. Along a subsequence, $Kz_j$ converges;
since $Tz_j=Tx_j\to y$, $z_j=Tz_j+Kz_j$ also converges to some $z$.
The continuity of $T$ gives $Tz=y$. Hence $\operatorname{Ran}T$ is
closed.

Let us prove the equivalence between injectivity and surjectivity. For
$j\geq0$, write
\[
 R_j:=\operatorname{Ran}T^j,
 \qquad N_j:=\ker T^j,
 \qquad R_0=X,\quad N_0=\{0\}.
\]
Each $N_j$ is closed. For $j\geq1$, the operator
\[
 L_j:=I-T^j
 =\sum_{l=1}^{j}(-1)^{l+1}\binom jl K^l
\]
is compact by Proposition~\ref{prop:propiedades-operadores-compactos}.
The closed range argument just proved, applied to
$T^j=I-L_j$, shows that each $R_j$ is closed. Moreover, $K$ commutes with
$T$ and therefore preserves all the spaces $R_j$ and $N_j$.

Suppose that $T$ is injective but not surjective. The inclusions
$R_{j+1}\subseteq R_j$ are then strict for every $j\geq0$.
Indeed, if $R_j=R_{j+1}$, then for each $x\in X$ there would exist $y\in X$ with
$T^jx=T^{j+1}y$; the injectivity of $T^j$ would imply $x=Ty$, a contradiction.
The separation fact gives vectors
\[
 x_j\in R_{j-1},\qquad \|x_j\|_X=1,
 \qquad\operatorname{dist}_X(x_j,R_j)>\frac{1}{2},
 \qquad j\geq1.
\]
If $k>j$, we have $Tx_j\in R_j$ and $Kx_k\in R_{k-1}\subseteq R_j$.
Consequently,
\[
 \|Kx_j-Kx_k\|_X
 =\|x_j-(Tx_j+Kx_k)\|_X>\frac{1}{2},
\]
contradicting the compactness of $K$. Thus injectivity implies
surjectivity.

Conversely, suppose that $T$ is surjective and $N_1\ne\{0\}$.
The inclusions $N_j\subseteq N_{j+1}$ are strict for every $j\geq0$.
If $N_j=N_{j+1}$ and $0\ne z\in N_1$, the surjectivity of $T^j$ would give
an $y$ with $T^jy=z$. Then $y\in N_{j+1}=N_j$, which would force
$z=0$, a contradiction. By separation, choose
\[
 x_j\in N_j,\qquad\|x_j\|_X=1,
 \qquad\operatorname{dist}_X(x_j,N_{j-1})>\frac{1}{2},
 \qquad j\geq1.
\]
If $j>k$, then $Tx_j\in N_{j-1}$ and
$Kx_k\in N_k\subseteq N_{j-1}$. The same identity gives
$\|Kx_j-Kx_k\|_X>\displaystyle\frac{1}{2}$, again contradicting compactness. This
proves the stated equivalence for every compact perturbation of
the identity.

It remains to compare the dimensions. By
Theorem~\ref{teo:schauder-adjunto-compacto}, $K'$ is compact on the Banach
space $X'$. The finite-dimensional kernel argument at the beginning of this proof,
applied to $I-K'=T'$, gives
\[
 m:=\dim\ker T'<\infty.
\]
Set $R:=\operatorname{Ran}T$ and $Z:=X/R$, with the quotient norm. The
map
\[
 Z'\longrightarrow\ker T',
 \qquad \ell\longmapsto\ell\circ\pi,
\]
where $\pi\colon X\to X/R$ is the projection, is a linear isomorphism.
Indeed, $T'\varphi=0$ is equivalent to $\varphi$ vanishing on $R$; in
that case, $\ell([x]):=\varphi(x)$ is well-defined and
\[
 |\ell([x])|
 =|\varphi(x-r)|
 \leq\|\varphi\|_{X'}\|x-r\|_X,
 \qquad r\in R.
\]
Taking the infimum gives continuity in the quotient norm. The
canonical injection $Z\to Z''$ is injective by
Corollary~\ref{cor:caracterizacion-dual-norma}. Since $\dim Z'=m$, it follows
that $Z$ is finite-dimensional; in finite dimensions, the dual has the
same dimension. Therefore,
\[
 \dim(X/R)=\dim Z=\dim Z'=m.
\]

Let $n:=\dim N$ and choose a basis $e_1,\ldots,e_n$ of $N$. The
coordinate functionals of this basis extend to functionals
$f_1,\ldots,f_n\in X'$ by
Corollary~\ref{cor:extension-hahn-banach-preserva-norma}. Then
\[
 Px:=\sum_{i=1}^{n}f_i(x)e_i,
 \qquad X_0:=\ker P,
\]
define a continuous projection onto $N$ and a decomposition
$X=N\oplus X_0$. The operator $T\restriction_{X_0}$ is a bijection from
$X_0$ onto $R$: it is injective because $N\cap X_0=\{0\}$ and
surjective because $T(x-Px)=Tx$.

Choose $y_1,\ldots,y_m\in X$ whose classes form a basis of $X/R$.
With $Y_0:=\operatorname{span}\{y_1,\ldots,y_m\}$, we have
$X=R\oplus Y_0$. Set $r:=\min\{n,m\}$ and define the continuous
finite-rank operator
\[
 Fx:=\sum_{i=1}^{r}f_i(x)y_i.
\]
Empty sums are interpreted as zero. By
Propositions~\ref{prop:b2-espacios-normados-operador-lineal-compacto-operador-acotado-decir}
and~\ref{prop:propiedades-operadores-compactos}, $K-F$ is compact. The
operator $S:=T+F=I-(K-F)$ satisfies the equivalence between injectivity and
surjectivity already proved.
For $\displaystyle x=x_0+\displaystyle\sum_{i=1}^{n}a_i e_i$, with $x_0\in X_0$, we have
\[
 Sx=Tx_0+\sum_{i=1}^{r}a_i y_i.
\]
The direct sum $X=R\oplus Y_0$ and the bijectivity of
$T\restriction_{X_0}$ give
\[
 \ker S=\operatorname{span}\{e_{r+1},\ldots,e_n\},
 \qquad
 \operatorname{Ran}S=R\oplus\operatorname{span}\{y_1,\ldots,y_r\}.
\]
If $n\leq m$, then $r=n$ and $S$ is injective; hence it is
surjective, which forces $r=m$. If $m\leq n$, then $r=m$ and $S$
is surjective; hence it is injective, which forces $r=n$. In both
cases, $n=m$, that is,
\[
 \dim\ker(I-K)=\dim\ker(I-K').
\]
This completes the proof.
\end{proof}

\begin{definition}\label{def:operador-fredholm-indice}\index{Fredholm operator}\index{Fredholm index}
Let $X$ and $Y$ be Banach spaces. An operator $T\in\mathcal L(X,Y)$ is called \textbf{Fredholm} if $\ker T$ is finite-dimensional, $\operatorname{Ran}(T)$ is closed, and $Y/\operatorname{Ran}(T)$ is finite-dimensional. Its \textbf{index} is
\[
\operatorname{ind}(T):=\dim\ker T-\dim\bigl(Y/\operatorname{Ran}(T)\bigr).
\]
In particular, Theorem~\ref{teo:riesz-schauder} implies that $I-K$ is Fredholm of index zero for every $K\in\mathcal K(X)$.
\end{definition}

\section{Hilbert spaces}\label{ap:espacios-de-hilbert}

Hilbert spaces combine completeness with a Euclidean geometry given by an inner product. This structure allows us to define orthogonal complements, projections, adjoints, and gradients, and provides the natural setting for many linear and quadratic variational problems. We mainly follow \cite[Section 1.2]{DrabekMilota2013}, \cite{brezis_functional_2011,kesavan_functional_2023,reed_simon_vol1}.

\begin{definition}\label{def:b3-espacios-hilbert-producto-interno}\index{inner product}
Let $V$ be a vector space over $\mathbb K=\mathbb R$ or $\mathbb C$. An \textbf{inner product} on $V$ is a map $\langle\cdot,\cdot\rangle_V\colon V\times V\longrightarrow\mathbb K$ such that, for all $x,y,z\in V$ and $\lambda\in\mathbb K$,
\begin{enumerate}[label=(\alph*)]
\item $\langle x+y,z\rangle_V=\langle x,z\rangle_V+\langle y,z\rangle_V$ and $\langle\lambda x,z\rangle_V=\lambda\langle x,z\rangle_V$;
\item $\langle x,y\rangle_V=\overline{\langle y,x\rangle_V}$;
\item $\langle x,x\rangle_V\geq0$, and $\langle x,x\rangle_V=0$ if and only if $x=0$.
\end{enumerate}
Thus, in the complex case, we adopt the convention of linearity in the first variable and antilinearity in the second.
\end{definition}

\begin{theorem}[Cauchy--Schwarz inequality]\label{teo:cauchy-schwarz-hilbert}\index{Cauchy--Schwarz inequality}
For any $x,y\in V$, we have $|\langle x,y\rangle_V|\leq\langle x,x\rangle_V^{\displaystyle\frac{1}{2}}\langle y,y\rangle_V^{\displaystyle\frac{1}{2}}$.
\end{theorem}

\begin{proof}
If $y=0$, the assertion is immediate. If $y\neq0$, positivity of the inner product applied to $x-\displaystyle\frac{\langle x,y\rangle_V}{\langle y,y\rangle_V}y$ gives
\[
0\leq\langle x,x\rangle_V-\frac{|\langle x,y\rangle_V|^2}{\langle y,y\rangle_V},
\]
from which the inequality follows.
\end{proof}

\begin{proposition}\label{prop:norma-inducida-producto-interno}
The formula $\|x\|_V:=\langle x,x\rangle_V^{\displaystyle\frac{1}{2}}$ defines a norm on $V$. Moreover, the parallelogram identity holds:
\[
\|x+y\|_V^2+\|x-y\|_V^2=2\|x\|_V^2+2\|y\|_V^2.
\]
\end{proposition}

\begin{proof}
Only the triangle inequality requires justification. By Cauchy--Schwarz,
\[
\|x+y\|_V^2
=\|x\|_V^2+2\operatorname{Re}\langle x,y\rangle_V+\|y\|_V^2
\leq(\|x\|_V+\|y\|_V)^2.
\]
The parallelogram identity follows by expanding both sides.
\end{proof}

\begin{definition}\label{def:b3-espacios-hilbert-espacio-hilbert-espacio-vectorial-producto-interno}\index{Hilbert space}
A vector space with an inner product is called a \textbf{Hilbert space} if it is complete in the norm induced by that inner product.
\end{definition}

\subsection{Orthogonality and projections}

\begin{definition}\label{def:complemento-ortogonal}\index{orthogonal complement}
Let $H$ be a Hilbert space and let $M\subseteq H$. The \textbf{orthogonal complement} of $M$ is
\[
M^{\perp}:=\{x\in H\mid \langle m,x\rangle_H=0\text{ for every }m\in M\}.
\]
Two subsets $A,B\subseteq H$ are orthogonal if $\langle a,b\rangle_H=0$ for all $a\in A$ and $b\in B$.
\end{definition}

\begin{proposition}\label{prop:propiedades-complemento-ortogonal}
For every $M\subseteq H$, the set $M^\perp$ is a closed subspace. If $M$ is a subspace, then $M\cap M^\perp=\{0\}$.
\end{proposition}

\begin{proof}
For each $m\in M$, the map $x\mapsto\langle m,x\rangle_H$ is continuous, so $M^\perp$ is an intersection of closed kernels. If $x\in M\cap M^\perp$, then $\|x\|_H^2=\langle x,x\rangle_H=0$, and therefore $x=0$.
\end{proof}

\begin{theorem}[Projection onto closed convex sets]\label{teo:proyeccion-convexo-cerrado-hilbert}\index{projection theorem}
Let $H$ be a Hilbert space and let $C\subseteq H$ be nonempty, closed, and convex. For each $x\in H$, there exists a unique $P_Cx\in C$ such that
\[
\|x-P_Cx\|_H=\inf_{z\in C}\|x-z\|_H.
\]
This element is characterized by
\[
\operatorname{Re}\langle x-P_Cx,z-P_Cx\rangle_H\leq0
\qquad
\text{for every }z\in C.
\]
\end{theorem}

\begin{proof}
Let $d:=\displaystyle\inf_{z\in C}\|x-z\|_H$ and choose $(z_n)$ in $C$ such that $\|x-z_n\|_H^2<d^2+\displaystyle\frac{1}{n}$. The parallelogram identity and convexity of $C$ imply
\[
\|z_n-z_m\|_H^2
=2\|x-z_n\|_H^2+2\|x-z_m\|_H^2-4\left\|x-\frac{z_n+z_m}{2}\right\|_H^2
\leq\frac{2}{n}+\frac{2}{m}.
\]
Therefore, $(z_n)$ is Cauchy and converges to some $y\in C$; necessarily, $\|x-y\|_H=d$. The same identity applied to two minimizers proves uniqueness.

Let $z\in C$. For $t\in(0,1)$, the point $y+t(z-y)$ belongs to $C$, so
\[
\|x-y\|_H^2
\leq\|x-y-t(z-y)\|_H^2
=\|x-y\|_H^2-2t\operatorname{Re}\langle x-y,z-y\rangle_H+t^2\|z-y\|_H^2.
\]
Dividing by $t$ and letting $t\to0^+$ gives the characterizing inequality. Conversely, this inequality implies $\|x-z\|_H^2\geq\|x-y\|_H^2$ for every $z\in C$.
\end{proof}

\begin{corollary}[Orthogonal decomposition]\label{cor:descomposicion-ortogonal-hilbert}\index{orthogonal decomposition}
If $M$ is a closed subspace of a Hilbert space $H$, then
\[
H=M\oplus M^\perp.
\]
For each $x\in H$, the element $P_Mx$ is the unique vector in $M$ such that $x-P_Mx\in M^\perp$. The operator $P_M\in\mathcal L(H)$ is called the \textbf{orthogonal projection} onto $M$ and satisfies $P_M^2=P_M$, $\|P_M\|=1$ if $M\neq\{0\}$, and
\[
\langle P_Mx,y\rangle_H=\langle x,P_My\rangle_H
\qquad
\text{for all }x,y\in H.
\]
\end{corollary}

\begin{proof}
Apply Theorem~\ref{teo:proyeccion-convexo-cerrado-hilbert} to $C=M$. Since $M$ is a subspace, the characterizing inequality applied to $z=P_Mx\pm m$ and, in the complex case, to $z=P_Mx\pm im$ implies $\langle m,x-P_Mx\rangle_H=0$ for every $m\in M$. Uniqueness of the decomposition is immediate. The linearity of $P_M$ follows from this characterization, and $\|P_Mx\|_H\leq\|x\|_H$ follows from the Pythagorean theorem. The last identity follows by decomposing $x$ and $y$ into their components in $M$ and $M^\perp$.
\end{proof}

\begin{corollary}\label{cor:doble-complemento-ortogonal}
For every subset $M\subseteq H$, we have
\[
(M^\perp)^\perp=\overline{\operatorname{span}(M)}.
\]
\end{corollary}

\begin{proof}
The inclusion $\overline{\operatorname{span}(M)}\subseteq(M^\perp)^\perp$ is immediate. Applying Corollary~\ref{cor:descomposicion-ortogonal-hilbert} to the closed subspace $N:=\overline{\operatorname{span}(M)}$, we write $x=P_Nx+(x-P_Nx)$. If $x\in(M^\perp)^\perp$, then $x-P_Nx\in N^\perp=M^\perp$, and therefore $\|x-P_Nx\|_H^2=\langle x,x-P_Nx\rangle_H=0$. Thus, $x=P_Nx\in N$.
\end{proof}

\subsection{Orthonormal systems and Riesz representation}

\begin{definition}\label{def:sistema-base-ortonormal}\index{orthonormal system}\index{orthonormal basis}
A family $(e_j)_{j\in J}$ in a Hilbert space $H$ is \textbf{orthonormal} if $\langle e_j,e_k\rangle_H=0$ for $j\neq k$ and $\|e_j\|_H=1$. It is called an \textbf{orthonormal basis} if $\overline{\operatorname{span}\{e_j\mid j\in J\}}=H$.
\end{definition}

\begin{theorem}[Bessel and Parseval]\label{teo:bessel-parseval}\index{Bessel's inequality}\index{Parseval's identity}
If $(e_j)_{j\in J}$ is an orthonormal family in $H$, then for every $x\in H$ and every finite subset $F\subseteq J$,
\[
\displaystyle\sum_{j\in F}|\langle x,e_j\rangle_H|^2\leq\|x\|_H^2.
\]
If $(e_j)_{j\in J}$ is an orthonormal basis, then
\[
x=\displaystyle\sum_{j\in J}\langle x,e_j\rangle_He_j,
\qquad
\|x\|_H^2=\displaystyle\sum_{j\in J}|\langle x,e_j\rangle_H|^2,
\]
where only countably many coefficients can be nonzero.
\end{theorem}

\begin{proof}
For finite $F\subseteq J$, put
$p_F:=\sum_{j\in F}\langle x,e_j\rangle_He_j$. Orthonormality gives
$\langle x-p_F,e_j\rangle_H=0$ for $j\in F$ and
\[
 \|x\|_H^2=\|x-p_F\|_H^2+
              \sum_{j\in F}|\langle x,e_j\rangle_H|^2.
\]
This proves Bessel's inequality. For every $m\geq1$, the set
\[
 J_m:=\{j\in J\mid |\langle x,e_j\rangle_H|\geq\frac{1}{m}\}
\]
is finite: by Bessel, every finite subset of $J_m$ has at most
$m^2\|x\|_H^2$ elements. Hence the set of nonzero coefficients,
$J_x=\bigcup_{m\geq1}J_m$, is countable.

If $J_x$ is infinite, enumerate it as $\{j_1,j_2,\ldots\}$.
Bessel shows that the series of squared coefficients converges,
and for $N>M$,
\[
 \left\|\sum_{\nu=M+1}^N\langle x,e_{j_\nu}\rangle_He_{j_\nu}
       \right\|_H^2
 =\sum_{\nu=M+1}^N|\langle x,e_{j_\nu}\rangle_H|^2.
\]
The partial sums are Cauchy and converge to a vector $p\in H$.
If $J_x$ is finite or empty, use its finite sum, or $p=0$.
In every case, continuity of the inner product gives
$\langle x-p,e_j\rangle_H=0$ for every $j\in J$.
Given $\varepsilon>0$ and a finite set $F_0\subseteq J_x$ with
$\sum_{j\in J_x\setminus F_0}|\langle x,e_j\rangle_H|^2<\varepsilon^2$,
the same tail identity gives $\|p-p_F\|_H<\varepsilon$ for every
finite $F\supseteq F_0$. This proves convergence of the finite sums
and independence of the order.

If the family is a basis, $x-p$ is orthogonal to a dense subspace
and hence, by continuity, to all of $H$. Thus $x-p=0$.
Passing to the limit in the Pythagorean identity yields
\[
 x=\sum_{j\in J}\langle x,e_j\rangle_He_j,
 \qquad
 \|x\|_H^2=\sum_{j\in J}|\langle x,e_j\rangle_H|^2.
\]
\end{proof}

\begin{theorem}[Riesz representation theorem]\label{teo:representacion-riesz-hilbert}\index{Riesz representation theorem}
Let $H$ be a Hilbert space and let $\ell\in H'$. There exists a unique $v\in H$ such that
\[
\ell(x)=\langle x,v\rangle_H
\qquad
\text{for every }x\in H.
\]
Moreover, $\|\ell\|_{H'}=\|v\|_H$. In the real case, the map $J\colon H\longrightarrow H'$, given by $J(v)(x)=\langle x,v\rangle_H$, is a linear isometry; in the complex case, it is an antilinear isometry.
\end{theorem}

\begin{proof}
If $\ell=0$, take $v=0$. Suppose that $\ell\neq0$ and let $M:=\ker\ell$, which is a closed subspace. By Corollary~\ref{cor:descomposicion-ortogonal-hilbert}, $H=M\oplus M^\perp$. The space $M^\perp$ is one-dimensional: the restriction of $\ell$ to $M^\perp$ is injective and takes values in $\mathbb K$. Choose $e\in M^\perp$ with $\|e\|_H=1$. Each $x\in H$ has a unique representation as $x=m+\alpha e$, with $m\in M$, and $\ell(x)=\alpha\ell(e)$.

In the real case, define $v:=\ell(e)e$, and in the complex case, $v:=\overline{\ell(e)}e$. With our convention of linearity in the first variable, $\langle x,v\rangle_H=\alpha\ell(e)=\ell(x)$. Uniqueness follows because $\langle x,v-w\rangle_H=0$ for every $x\in H$ implies $v=w$. Cauchy--Schwarz gives $\|\ell\|\leq\|v\|_H$, while, if $v\neq0$, evaluating at $\displaystyle\frac{v}{\|v\|_H}$ gives the reverse inequality.
\end{proof}

\begin{corollary}[Complex-linear form of Riesz representation]
\label{cor:riesz-lineal-desde-conjugado}
Let $H$ be a complex Hilbert space with inner product linear in the
first variable. Then
\[
 \mathscr R_H\colon\overline H\longrightarrow H',
 \qquad
 \mathscr R_H(\overline v)(x):=\langle x,v\rangle_H,
\]
is a complex-linear isometric isomorphism. In particular, every
complex-linear identification of a complex Hilbert space with its dual
obtained from its inner product must have the conjugate
space $\overline H$ as its domain; the direct map $H\to H'$ is antilinear.
\end{corollary}

\begin{proof}
Theorem~\ref{teo:representacion-riesz-hilbert} proves bijectivity and
equality of norms. If $\lambda\in\mathbb C$, scalar multiplication in
$\overline H$ satisfies $\lambda\overline v=\overline{\overline\lambda v}$;
therefore,
\[
 \mathscr R_H(\lambda\overline v)(x)
 =\langle x,\overline\lambda v\rangle_H
 =\lambda\langle x,v\rangle_H
 =\lambda\mathscr R_H(\overline v)(x).
\]
\end{proof}

\begin{theorem}[Lax--Milgram theorem]\label{teo:lax-milgram}\index{Lax--Milgram theorem}
Let $H$ be a real Hilbert space and let $a\colon H\times H\longrightarrow\mathbb R$ be bilinear. Suppose there exist constants $C\geq0$ and $\alpha>0$ such that
\[
|a(u,v)|\leq C\|u\|_H\|v\|_H,
\qquad
a(u,u)\geq\alpha\|u\|_H^2
\]
for all $u,v\in H$. Then, for every $F\in H'$, there exists a unique $u\in H$ such that $a(u,v)=F(v)$ for every $v\in H$. Moreover, $\|u\|_H\leq\displaystyle\frac{\|F\|_{H'}}{\alpha}$.
\end{theorem}

\begin{proof}
By Theorem~\ref{teo:representacion-riesz-hilbert}, for each $u\in H$ there exists a unique $Tu\in H$ such that $a(u,v)=\langle v,Tu\rangle_H$ for every $v\in H$. The operator $T$ is linear and continuous. Coercivity implies
\[
\alpha\|u\|_H^2\leq\langle u,Tu\rangle_H\leq\|u\|_H\|Tu\|_H,
\]
so $\|Tu\|_H\geq\alpha\|u\|_H$. In particular, $T$ is injective and its range is closed. If $w\perp\operatorname{Ran}(T)$, then $0=\langle w,Tw\rangle_H=a(w,w)\geq\alpha\|w\|_H^2$, so $w=0$. Therefore, the range of $T$ is dense and closed, hence is all of $H$.

Let $f\in H$ be the Riesz representative of $F$. There exists a unique $u\in H$ such that $Tu=f$, and then $a(u,v)=\langle v,f\rangle_H=F(v)$. Finally, $\alpha\|u\|_H\leq\|Tu\|_H=\|f\|_H=\|F\|_{H'}$.
\end{proof}

\subsection{Adjoint operators on Hilbert spaces}

\begin{theorem}\label{teo:adjunto-operador-acotado-hilbert}
Let $H_1$ and $H_2$ be Hilbert spaces and let $T\in\mathcal L(H_1,H_2)$. There exists a unique operator $T^*\in\mathcal L(H_2,H_1)$ such that
\[
\langle Tx,y\rangle_{H_2}=\langle x,T^*y\rangle_{H_1}
\qquad
\text{for all }x\in H_1,\,y\in H_2.
\]
Moreover, $\|T^*\|=\|T\|$, $(S\circ T)^*=T^*\circ S^*$, and $(T^*)^*=T$.
\end{theorem}

\begin{proof}
For each $y\in H_2$, the map $x\mapsto\langle Tx,y\rangle_{H_2}$ belongs to $H_1'$. Theorem~\ref{teo:representacion-riesz-hilbert} gives a unique vector $T^*y\in H_1$ satisfying the identity. Uniqueness allows us to verify the linearity of $T^*$, and Cauchy--Schwarz gives $\|T^*y\|_{H_1}\leq\|T\|\|y\|_{H_2}$. Applying the same construction to $T^*$ and using the defining identity gives $(T^*)^*=T$, and therefore $\|T\|\leq\|T^*\|$. The remaining identities are verified using the inner product.
\end{proof}

\begin{definition}\label{def:operadores-autoadjuntos-positivos-unitarios}\index{self-adjoint operator}\index{positive operator}\index{unitary operator}
Let $H$ be a Hilbert space and let $T\in\mathcal L(H)$.
\begin{enumerate}[label=(\alph*)]
\item $T$ is \textbf{self-adjoint} if $T=T^*$.
\item $T$ is \textbf{positive} if $\langle Tx,x\rangle_H\geq0$ for every $x\in H$.
\item $T$ is \textbf{unitary} if $T^*T=TT^*=I$.
\end{enumerate}
\end{definition}

For unbounded operators, the definition of the adjoint must include its domain.

\begin{definition}\label{def:adjunto-operador-densamente-definido}\index{adjoint operator!densely defined}
Let $H_1$ and $H_2$ be Hilbert spaces and let $A\colon \mathcal D(A)\subseteq H_1\longrightarrow H_2$ be a densely defined linear operator. Define
\[
\mathcal D(A^*):=\left\{y\in H_2\middle|\text{there exists }z\in H_1\text{ such that }\langle Ax,y\rangle_{H_2}=\langle x,z\rangle_{H_1}\text{ for every }x\in\mathcal D(A)\right\}.
\]
For $y\in\mathcal D(A^*)$, the vector $z$ is unique, and we define $A^*y:=z$.
\end{definition}

\begin{proposition}\label{prop:adjunto-siempre-cerrado}
The adjoint of every densely defined linear operator is closed.
\end{proposition}

\begin{proof}
Let $y_n\in\mathcal D(A^*)$ satisfy $y_n\to y$ in $H_2$ and $A^*y_n\to z$ in $H_1$. For each $x\in\mathcal D(A)$,
\[
\langle Ax,y\rangle_{H_2}
=\lim_{n\to\infty}\langle Ax,y_n\rangle_{H_2}
=\lim_{n\to\infty}\langle x,A^*y_n\rangle_{H_1}
=\langle x,z\rangle_{H_1}.
\]
Therefore, $y\in\mathcal D(A^*)$ and $A^*y=z$.
\end{proof}

\begin{theorem}\label{teo:cerrabilidad-mediante-adjunto}
Let $A\colon \mathcal D(A)\subseteq H_1\longrightarrow H_2$ be linear and densely defined. Then $A$ is closable if and only if $\mathcal D(A^*)$ is dense in $H_2$. In this case, $\overline A=A^{**}$.
\end{theorem}

The proof relies on identifying $H_j$ with $H_j'$ by Theorem~\ref{teo:representacion-riesz-hilbert} and on the relation between the graph of $A$ and its orthogonal complement. It can be found in \cite{reed_simon_vol1,brezis_functional_2011}.

\begin{definition}\label{def:operador-simetrico-autoadjunto-no-acotado}\index{symmetric operator}\index{self-adjoint operator}
Let $H$ be a Hilbert space and let $A\colon \mathcal D(A)\subseteq H\longrightarrow H$ be linear and densely defined. We say that $A$ is \textbf{symmetric} if $\langle Ax,y\rangle_H=\langle x,Ay\rangle_H$ for all $x,y\in\mathcal D(A)$, which is equivalent to $A\subseteq A^*$. We say that $A$ is \textbf{self-adjoint} if $\mathcal D(A)=\mathcal D(A^*)$ and $A=A^*$ on this domain.
\end{definition}

\begin{proposition}\label{prop:b3-espacios-hilbert-espacio-hilbert-operador-simetrico-siguientes-equivalentes}
Let $H$ be a complex Hilbert space and let $A\colon \mathcal D(A)\subseteq H\longrightarrow H$ be a densely defined symmetric operator. The following assertions are equivalent:
\begin{enumerate}[label=(\alph*)]
\item $A$ is self-adjoint.
\item $A$ is closed and $\ker(A^*-iI)=\ker(A^*+iI)=\{0\}$.
\item $\operatorname{Ran}(A-iI)=\operatorname{Ran}(A+iI)=H$.
\end{enumerate}
\end{proposition}

\begin{proof}
For every $x\in\mathcal D(A)$, symmetry implies
\[
\|(A\pm iI)x\|_H^2=\|Ax\|_H^2+\|x\|_H^2.
\]
If $A$ is closed, this estimate shows that the ranges of $A\pm iI$ are closed. Moreover,
\[
\operatorname{Ran}(A+iI)^\perp=\ker(A^*-iI),
\qquad
\operatorname{Ran}(A-iI)^\perp=\ker(A^*+iI).
\]
If $A$ is self-adjoint, it is closed and the preceding kernels are trivial; therefore, $(a)\Rightarrow(b)$. If $(b)$ holds, both ranges are closed and their orthogonal complements are trivial, so both equal $H$. Thus, $(b)\Rightarrow(c)$.

Suppose $(c)$ and let $y\in\mathcal D(A^*)$. Choose $x\in\mathcal D(A)$ such that $(A+iI)x=(A^*+iI)y$. Since $A\subseteq A^*$, we have $(A^*+iI)(y-x)=0$. Now $\ker(A^*+iI)=\operatorname{Ran}(A-iI)^\perp=\{0\}$, so $y=x\in\mathcal D(A)$. Consequently, $\mathcal D(A^*)\subseteq\mathcal D(A)$ and $A$ is self-adjoint. This proves $(c)\Rightarrow(a)$.
\end{proof}

\subsection{Compact self-adjoint operators}

\begin{proposition}\label{prop:compacto-hilbert-limite-rango-finito}
Let $H_1$ and $H_2$ be Hilbert spaces. An operator $K\in\mathcal L(H_1,H_2)$ is compact if and only if there exists a sequence of finite-rank operators $(K_n)$ such that $\|K_n-K\|_{\mathcal L(H_1,H_2)}\to0$.
\end{proposition}

\begin{proof}
If $K$ is the norm limit of finite-rank operators, then it is compact by Propositions~\ref{prop:b2-espacios-normados-operador-lineal-compacto-operador-acotado-decir} and \ref{prop:propiedades-operadores-compactos}.

Conversely, suppose that $K$ is compact and set
\[
C:=\overline{K(B_{d_{H_1}}(0,1))},
\qquad
M:=\overline{\operatorname{span}(C)}.
\]
The set $C$ is compact because $K(B_{d_{H_1}}(0,1))$ is relatively compact. For each $r\in\mathbb N$, choose a finite set $F_r\subseteq C$ such that the balls of radius $\displaystyle\frac{1}{r}$ centered at the points of $F_r$ cover $C$. The union $D:=\displaystyle\bigcup_{r=1}^{\infty}F_r$ is countable and dense in $C$.

Enumerate $D=\{y_1,y_2,\dots\}$ and let $\mathbb K_0:=\mathbb Q$ in the real case and $\mathbb K_0:=\mathbb Q+i\mathbb Q$ in the complex case. The set
\[
D_0
:=
\left\{
\displaystyle\sum_{j=1}^{s}a_jy_{n_j}
\middle|
 s\in\mathbb N,\ n_1,\dots,n_s\in\mathbb N,\ a_1,\dots,a_s\in\mathbb K_0
\right\}
\]
is countable, since it is a countable union of images of the countable sets $\mathbb N^s\times\mathbb K_0^s$.

Let us verify that $D_0$ is dense in $M$. It suffices to approximate the elements of $\operatorname{span}(C)$, since this subspace is dense in $M$ by definition. Let $z=\displaystyle\sum_{j=1}^{s}\alpha_jc_j$, with $c_j\in C$, $\alpha_j\in\mathbb K$, and let $\varepsilon>0$. For each $j$, choose $y_{n_j}\in D$ sufficiently close to $c_j$ that
\[
\displaystyle\sum_{j=1}^{s}|\alpha_j|\,\|c_j-y_{n_j}\|_{H_2}<\frac{\varepsilon}{2}.
\]
Then choose $a_j\in\mathbb K_0$ so that
\[
\displaystyle\sum_{j=1}^{s}|\alpha_j-a_j|\,\|y_{n_j}\|_{H_2}<\frac{\varepsilon}{2};
\]
if any of the vectors $y_{n_j}$ is zero, the corresponding term imposes no condition on $a_j$. Then $w:=\displaystyle\sum_{j=1}^{s}a_jy_{n_j}$ belongs to $D_0$ and
\[
\|z-w\|_{H_2}
\leq
\displaystyle\sum_{j=1}^{s}|\alpha_j|\,\|c_j-y_{n_j}\|_{H_2}
+
\displaystyle\sum_{j=1}^{s}|\alpha_j-a_j|\,\|y_{n_j}\|_{H_2}
<\varepsilon.
\]
Consequently, $D_0$ is countable and dense in $M$, so $M$ is separable.

For every $x\in H_1$, we have $Kx\in M$: if $x\neq0$, simply write $Kx=2\|x\|_{H_1}K\left(\displaystyle\frac{x}{2\|x\|_{H_1}}\right)$, and the vector in parentheses belongs to $B_{d_{H_1}}(0,1)$; the case $x=0$ is immediate. If $M$ is finite-dimensional, this shows that $K$ already has finite rank. Thus, suppose that $M$ is infinite-dimensional. Applying the Gram--Schmidt procedure to a dense sequence in $M$, and at each step discarding the vectors that belong to the subspace spanned by the preceding ones, gives a countable orthonormal basis $(e_j)_{j\in\mathbb N}$ of $M$. Let $P_n$ be the orthogonal projection from $H_2$ onto $\operatorname{span}\{e_1,\dots,e_n\}$. Since $K(H_1)\subseteq M$, the operator $K_n:=P_nK$ has finite rank.

By Theorem~\ref{teo:bessel-parseval}, $P_ny\to y$ for every $y\in M$. This convergence is uniform on $C$. To verify this, let $\varepsilon>0$ and choose a finite net $\{z_1,\dots,z_s\}\subseteq C$ of radius $\displaystyle\frac{\varepsilon}{2}$. There exists $N\in\mathbb N$ such that $\|(I-P_n)z_i\|_{H_2}<\displaystyle\frac{\varepsilon}{2}$ for every $i\in\{1,\dots,s\}$ and $n\geq N$. Given $y\in C$, choose $z_i$ with $\|y-z_i\|_{H_2}<\displaystyle\frac{\varepsilon}{2}$. Since $I-P_n$ is an orthogonal projection and therefore has norm at most one,
\[
\|(I-P_n)y\|_{H_2}
\leq
\|(I-P_n)(y-z_i)\|_{H_2}
+
\|(I-P_n)z_i\|_{H_2}
<\varepsilon.
\]
Consequently,
\[
\|K-K_n\|
=
\sup_{\substack{x\in H_1\\\|x\|_{H_1}\leq1}}\|(I-P_n)Kx\|_{H_2}
\leq
\sup_{y\in C}\|(I-P_n)y\|_{H_2}
\longrightarrow0.
\]
\end{proof}

\begin{proposition}[Bilinear forms and self-adjoint operators]
\label{prop:formas-bilineales-simetricas-hilbert}
\index{bilinear form!representation by an operator}
\index{self-adjoint operator!representation of bilinear forms}
Let $H$ be a real Hilbert space. For every continuous bilinear form $B\colon H\times H\longrightarrow\mathbb R$, there exists a unique operator $A_B\in\mathcal L(H)$ such that
\[
B(x,y)=\langle A_Bx,y\rangle_H
\]
for all $x,y\in H$. The correspondence $B\longmapsto A_B$ is a linear isometric isomorphism between $\mathcal L^2(H;\mathbb R)$ and $\mathcal L(H)$. Moreover, $B$ is symmetric if and only if $A_B$ is self-adjoint.

Define the rank of $B$ to be the rank of $A_B$. With this convention, symmetric bilinear forms of finite rank correspond exactly to self-adjoint operators of finite rank. If $\mathcal F_s(H)$ denotes the subspace of these forms, then its closure in $\mathcal L_s^2(H;\mathbb R)$ corresponds to the space of compact self-adjoint operators. In particular, if $H$ is infinite-dimensional and $B_0(x,y):=\langle x,y\rangle_H$, then $B_0\notin\overline{\mathcal F_s(H)}$.
\end{proposition}

\begin{proof}
For fixed $x\in H$, the map $y\longmapsto B(x,y)$ belongs to $H'$. Theorem~\ref{teo:representacion-riesz-hilbert} gives a unique vector $A_Bx\in H$ satisfying $B(x,y)=\langle A_Bx,y\rangle_H$ for every $y\in H$. Uniqueness of the Riesz representatives shows that $A_B$ is linear. Similarly,
\[
\|A_Bx\|_H
=
\sup_{\substack{y\in H\\\|y\|_H\leq1}}|B(x,y)|
\leq
\|B\|\,\|x\|_H,
\]
so $A_B\in\mathcal L(H)$ and $\|A_B\|\leq\|B\|$. The reverse inequality follows from $|B(x,y)|\leq\|A_B\|\|x\|_H\|y\|_H$, so $\|A_B\|=\|B\|$. Linearity and surjectivity of the correspondence follow by associating the form $B_A(x,y):=\langle Ax,y\rangle_H$ to $A\in\mathcal L(H)$.

Symmetry of $B$ is equivalent to the identity $\langle A_Bx,y\rangle_H=\langle x,A_By\rangle_H$ for all $x,y\in H$, which is precisely condition $A_B=A_B^*$ of Theorem~\ref{teo:adjunto-operador-acotado-hilbert}. The assertion about rank is then immediate from the definition adopted.

It remains to identify the closure of $\mathcal F_s(H)$. Every norm limit of finite-rank self-adjoint operators is compact by Proposition~\ref{prop:propiedades-operadores-compactos} and self-adjoint by the continuity of the adjoint established in Theorem~\ref{teo:adjunto-operador-acotado-hilbert}. Conversely, let $K\in\mathcal K(H)$ be self-adjoint. Proposition~\ref{prop:compacto-hilbert-limite-rango-finito} gives finite-rank operators $T_n$ such that $T_n\to K$ in norm. The adjoint $T_n^*$ also has finite rank. Indeed, $\operatorname{Ran}(T_n)$ is closed because it is finite-dimensional; if $P_n$ is the orthogonal projection onto this range, then $T_n=P_nT_n$, and therefore $T_n^*=T_n^*P_n$. Thus the range of $T_n^*$ is contained in the image under $T_n^*$ of the finite-dimensional space $\operatorname{Ran}(P_n)$. Consequently, the operators $S_n:=\displaystyle\frac{1}{2}(T_n+T_n^*)$ are self-adjoint and have finite rank. Moreover, using $K^*=K$ and $\|T_n^*-K^*\|=\|T_n-K\|$ gives $\|S_n-K\|\leq\|T_n-K\|\to0$.

Finally, if $H$ is infinite-dimensional, an orthonormal sequence $(e_n)_{n\in\mathbb N}$ can be constructed inductively. Since $\|e_n-e_m\|_H=\sqrt2$ for $n\neq m$, the sequence $(e_n)$ has no convergent subsequences. The identity on $H$ is therefore not compact. Since $B_0$ corresponds to the identity, we conclude that $B_0\notin\overline{\mathcal F_s(H)}$.
\end{proof}

\begin{corollary}\label{cor:forma-simetrizada-rango-finito-hilbert}
Let $H$ be a real Hilbert space and let $\ell,m\in H'$. The symmetric bilinear form
\[
B_{\ell,m}:=\ell\otimes m+m\otimes\ell,
\qquad
B_{\ell,m}(x,y)=\ell(x)m(y)+m(x)\ell(y),
\]
has rank at most two.
\end{corollary}

\begin{proof}
Let $v_\ell,v_m\in H$ be the Riesz representatives of $\ell$ and $m$. The operator associated with $B_{\ell,m}$ is $x\longmapsto\ell(x)v_m+m(x)v_\ell$, whose image is contained in $\operatorname{span}\{v_\ell,v_m\}$.
\end{proof}

\begin{theorem}[Spectral theorem for compact self-adjoint operators]\label{teo:espectral-compacto-autoadjunto}\index{spectral theorem!compact self-adjoint operators}
Let $H$ be a Hilbert space and let $K\in\mathcal K(H)$ be self-adjoint. The nonzero eigenvalues of $K$ are real, have finite multiplicity, and can accumulate only at $0$. There exists an orthonormal family $(e_j)_{j\in J}$ of eigenvectors, with real eigenvalues $(\lambda_j)_{j\in J}$, such that
\[
\overline{\operatorname{Ran}(K)}=\overline{\operatorname{span}\{e_j\mid j\in J\}},
\qquad
Kx=\displaystyle\sum_{j\in J}\lambda_j\langle x,e_j\rangle_He_j
\]
for every $x\in H$, and for each $\varepsilon>0$ there are only finitely many indices $j$ with $|\lambda_j|\geq\varepsilon$. In particular, if $H$ is separable, this family can be extended by an orthonormal basis of $\ker K$ to obtain an orthonormal basis of $H$ consisting of eigenvectors of $K$.
\end{theorem}

\begin{proof}
For a self-adjoint operator, $q(x):=\langle Kx,x\rangle_H$ is real.
Set $\alpha:=\sup_{\|x\|_H=1}|q(x)|$. Then $\alpha=\|K\|$.
One inequality follows from Cauchy--Schwarz.
For the other, given unit vectors $x,y$, multiply $y$ by a scalar
of modulus one, if necessary, so that $\langle Kx,y\rangle_H$
is real and nonnegative. Real polarization gives
\[
 |\langle Kx,y\rangle_H|
 =\frac{1}{4}\bigl(q(x+y)-q(x-y)\bigr)
 \leq\frac{\alpha}{4}
          \bigl(\|x+y\|_H^2+\|x-y\|_H^2\bigr)=\alpha.
\]
Taking the supremum in $y$ and then in $x$ yields $\|K\|\leq\alpha$.

The case $K=0$ is immediate. Suppose $K\neq0$. There is a sequence of unit vectors $x_j$ and,
after taking a subsequence, a sign such that $q(x_j)\to\lambda$,
where $\lambda=\|K\|$ or $\lambda=-\|K\|$.
Since $\lambda$ is real,
\[
 \|Kx_j-\lambda x_j\|_H^2
 =\|Kx_j\|_H^2-2\lambda q(x_j)+\lambda^2
 \leq2\|K\|^2-2\lambda q(x_j)\longrightarrow0.
\]
By compactness, we may assume $Kx_j\to z$ in norm.
Since $\lambda\neq0$, also $x_j\to e:=\displaystyle\frac{z}{\lambda}$.
Thus $\|e\|_H=1$ and $Ke=\lambda e$: the operator norm is
attained as the absolute value of an eigenvalue.

Every eigenvalue of a self-adjoint operator is real:
$Ke=\mu e$ gives
$\mu\|e\|_H^2=\langle Ke,e\rangle_H\in\mathbb R$.
If $Ke=\mu e$ and $Kf=\nu f$, then
$(\mu-\nu)\langle e,f\rangle_H=0$; distinct eigenspaces are
orthogonal. An eigenspace with $\mu\neq0$ cannot be
infinite dimensional: an orthonormal sequence in it would have
images separated by $|\mu|\sqrt2$, contradicting compactness.

Choose $e_1$ as above and restrict $K$ to $e_1^\perp$.
This subspace is invariant since
$\langle Kx,e_1\rangle_H=\langle x,Ke_1\rangle_H=0$ for
$x\perp e_1$. The restriction remains compact and self-adjoint.
If it is nonzero, choose $e_2$ there with $|\lambda_2|$ equal
to its norm, and continue on
$H_N:=\operatorname{span}\{e_1,\ldots,e_N\}^\perp$.
The procedure terminates with a zero restriction or gives
an orthonormal sequence with nonincreasing $|\lambda_j|$.
In the latter case $|\lambda_j|\to0$: otherwise the images
$Ke_j=\lambda_je_j$ would contain a subsequence separated
by a positive constant.

Let $S$ be the closed span of the chosen vectors.
Every kernel element is orthogonal to them. Conversely, if
$x\perp S$, in the infinite case
\[
 \|Kx\|_H\leq\|K|_{H_N}\|\,\|x\|_H
             =|\lambda_{N+1}|\|x\|_H\longrightarrow0;
\]
in the finite case the remaining restriction is already zero.
Hence $S^\perp=\ker K$. Decomposing $x$ in $S\oplus\ker K$
and using continuity gives
\[
 Kx=\sum_j\lambda_j\langle x,e_j\rangle_He_j
\]
with norm convergence. Each $e_j$ belongs to the range of $K$
since $e_j=K(\displaystyle\frac{e_j}{\lambda_j})$, and the series shows that the
range is contained in $S$. Thus
$\overline{\operatorname{Ran}K}=S$.
For every $\varepsilon>0$, only finitely many terms satisfy
$|\lambda_j|\geq\varepsilon$.
The complete description $H=S\oplus\ker K$ excludes any other
nonzero eigenvalue. Completing with an orthonormal basis of
the kernel proves the final assertion.
\end{proof}

Consultation sources: \cite{DrabekMilota2013,reed_simon_vol1,brezis_functional_2011}.

\subsection{Weak and strong regularity of Hilbert-space-valued maps}

Scalar regularity of all the components of a Hilbert-space-valued map determines almost all of its norm regularity. Although the conclusion is formulated in terms of differentiability, the argument relies on the uniform boundedness principle and Riesz representation; it does not require differential calculus on Banach manifolds. For this reason, we include the result in this section.

\begin{theorem}[Weak and strong differentiability]\label{teo:diferenciabilidad-debil-fuerte-hilbert}
Let $U\subseteq\mathbb R^n$ be open, let $H$ be a Hilbert space, and let $k\in\mathbb N$. Suppose that $F\colon U\longrightarrow H$ is of class $C^k$ in the weak sense, that is, for each $h\in H$ the scalar function $x\mapsto\langle F(x),h\rangle_H$ belongs to $C^k(U)$. Then $F$ is of class $C^{k-1}$ in norm. More precisely, for each multiindex $\alpha$ with $|\alpha|\leq k$ there exists a unique map $F_\alpha\colon U\longrightarrow H$ such that
\[
\langle F_\alpha(x),h\rangle_H
=
D^\alpha\langle F(\,\cdot\,),h\rangle_H(x),
\qquad h\in H,
\]
and, if $|\alpha|\leq k-1$, the map $F_\alpha$ is norm-continuous and equals $D^\alpha F$. In particular, every weakly smooth map with values in $H$ is smooth in norm. The same assertion holds for maps defined on a finite-dimensional manifold with or without boundary, after working in charts.
\end{theorem}

Another formulation of this result can be found in
\cite[Appendix B.5]{Guneysu2017Covariant}.

\begin{proof}
For $h\in H$, write $\phi_h(y):=\langle F(y),h\rangle_H$. We begin by constructing the vectors $F_\alpha(x)$. If $\alpha=0$, set $F_0:=F$. Now fix $x\in U$ and a multiindex $\alpha$ of length $r:=|\alpha|\geq1$. Choose indices $i_1,\dots,i_r\in\{1,\dots,n\}$ in which each $i$ occurs exactly $\alpha_i$ times and, for a function $G$ defined near $y$, denote its finite difference in the direction $e_i$ by $\Delta_{t,i}G(y):=G(y+te_i)-G(y)$.

Take a sequence $t_\nu>0$ tending to zero and small enough that all the points appearing in the following differences belong to $U$. Define
\[
v_\nu
:=
\frac{1}{t_\nu^r}
\Delta_{t_\nu,i_1}\cdots\Delta_{t_\nu,i_r}F(x)
\in H.
\]
The fundamental theorem of calculus, applied successively to the scalar function $\phi_h$, gives
\[
\langle v_\nu,h\rangle_H
=
\int_{[0,1]^r}
D^\alpha\phi_h\left(x+t_\nu\displaystyle\sum_{j=1}^rs_je_{i_j}\right)
\,ds_1\cdots ds_r.
\]
Since $D^\alpha\phi_h$ is continuous, it follows that $\langle v_\nu,h\rangle_H\longrightarrow D^\alpha\phi_h(x)$ for every $h\in H$. In particular, the sequence of continuous antilinear functionals $T_\nu\colon H\longrightarrow\mathbb K$, given by $T_\nu(h):=\langle v_\nu,h\rangle_H$, is pointwise bounded. The Banach--Steinhaus theorem, Theorem~\ref{teo: banach steinhaus}, applied to the linear functionals $h\mapsto\overline{T_\nu(h)}$, implies that there exists $C_{x,\alpha}>0$ such that $\|T_\nu\|\leq C_{x,\alpha}$ for every $\nu$. Consequently, the pointwise limit
\[
T_{x,\alpha}(h):=D^\alpha\phi_h(x)
\]
is antilinear and continuous and satisfies $|T_{x,\alpha}(h)|\leq C_{x,\alpha}\|h\|_H$. The antilinear version of the Riesz representation theorem, Theorem~\ref{teo:representacion-riesz-hilbert}, obtained by conjugating the functional, determines a unique vector $F_\alpha(x)\in H$ such that $T_{x,\alpha}(h)=\langle F_\alpha(x),h\rangle_H$ for every $h\in H$. Since $x$ was arbitrary, this constructs $F_\alpha\colon U\longrightarrow H$ for each $|\alpha|\leq k$. Uniqueness follows because the inner product separates points. Moreover, the identity
\[
\langle F_\alpha(x),h\rangle_H=D^\alpha\phi_h(x)
\]
shows that each $F_\alpha$ is weakly continuous.

We now prove norm continuity for $|\alpha|\leq k-1$. Fix a closed convex Euclidean ball $K\subseteq U$. For each $h\in H$ and each $i\in\{1,\dots,n\}$, the function $z\mapsto\langle F_{\alpha+e_i}(z),h\rangle_H=D^{\alpha+e_i}\phi_h(z)$ is continuous and therefore bounded on $K$. Applying Banach--Steinhaus to the linear functionals $h\mapsto\overline{\langle F_{\alpha+e_i}(z),h\rangle_H}$, with $z\in K$ and $i\in\{1,\dots,n\}$, gives a constant $C_{K,\alpha}>0$ such that
\[
\|F_{\alpha+e_i}(z)\|_H
\leq
C_{K,\alpha}
\]
for every $z\in K$ and every $i$. If $x,y\in K$, the segment $x+t(y-x)$ is contained in $K$. The fundamental theorem of calculus applied to $D^\alpha\phi_h$ gives
\[
\langle F_\alpha(y)-F_\alpha(x),h\rangle_H
=
\int_0^1
\displaystyle\sum_{i=1}^n(y_i-x_i)
\langle F_{\alpha+e_i}(x+t(y-x)),h\rangle_H\,dt.
\]
Taking the supremum over $\|h\|_H\leq1$ yields
\[
\|F_\alpha(y)-F_\alpha(x)\|_H
\leq
C_{K,\alpha}\displaystyle\sum_{i=1}^n|y_i-x_i|
\leq
\sqrt n\,C_{K,\alpha}\|y-x\|.
\]
Thus, $F_\alpha$ is locally Lipschitz and, in particular, norm-continuous for $|\alpha|\leq k-1$.

It remains to identify these maps with the strong derivatives of $F$. Let $|\alpha|\leq k-2$ and define $A_\alpha(x)\in\mathcal L(\mathbb R^n,H)$ by
\[
A_\alpha(x)v
:=
\displaystyle\sum_{i=1}^nv_iF_{\alpha+e_i}(x).
\]
Fix $x\in U$ and take $v$ small enough that the segment $x+tv$, with $0\leq t\leq1$, is contained in a closed ball $K\subseteq U$. For every $h\in H$, the fundamental theorem of calculus gives
\[
\begin{split}
&\left\langle F_\alpha(x+v)-F_\alpha(x)-A_\alpha(x)v,h\right\rangle_H\\
&\qquad=
\int_0^1\displaystyle\sum_{i=1}^nv_i
\left\langle F_{\alpha+e_i}(x+tv)-F_{\alpha+e_i}(x),h\right\rangle_H\,dt.
\end{split}
\]
Taking the supremum over $\|h\|_H\leq1$ again, we obtain
\[
\frac{\|F_\alpha(x+v)-F_\alpha(x)-A_\alpha(x)v\|_H}{\|v\|_{\mathbb R^n}}
\leq
\sqrt n\max_{1\leq i\leq n}\sup_{0\leq t\leq1}
\|F_{\alpha+e_i}(x+tv)-F_{\alpha+e_i}(x)\|_H.
\]
The right-hand side tends to zero as $v\to0$, because the maps $F_{\alpha+e_i}$ are norm-continuous. This proves that $F_\alpha$ is Fréchet differentiable and that $DF_\alpha(x)=A_\alpha(x)$. Since the $F_{\alpha+e_i}$ are continuous, $x\mapsto A_\alpha(x)$ is also continuous in the operator norm. Induction on $|\alpha|$ shows that $F\in C^{k-1}(U,H)$ and that its coordinate derivatives are the corresponding maps $F_\alpha$.

In interior charts the preceding argument applies. In a boundary chart we work in a half-space and choose convex closed half-balls. At boundary points, finite differences in the normal direction are taken inward; tangential differences remain in the chart for sufficiently small steps. Continuity up to the boundary of the scalar derivatives constructs the same vectors $F_\alpha$ there. Banach--Steinhaus gives the same bounds on each compact half-ball. The integral identities along segments, which remain in the half-space, prove norm continuity of $F_\alpha$ for $|\alpha|\leq k-1$ and existence of the lower-order derivatives with continuous values up to the boundary.

To specify regularity in a boundary chart, put $\ell=k-1$ and write its coordinates as $(x',t)$, with $t\geq0$. Choose numbers $a_1,\ldots,a_{\ell+1}$ satisfying
\[
\sum_{j=1}^{\ell+1}a_j(-j)^r=1,\qquad 0\leq r\leq\ell.
\]
They exist because the corresponding Vandermonde matrix is invertible. After shrinking the chart, define for $t<0$
\[
\widetilde F(x',t)=\sum_{j=1}^{\ell+1}a_jF(x',-jt),
\]
and retain $\widetilde F=F$ for $t\geq0$. The preceding identities match all derivatives of order at most $\ell$ at $t=0$. This gives a norm $C^\ell$ extension. Regularity is local and invariant under smooth coordinate changes, so the conclusions patch together on the manifold. If $F$ is weakly smooth, the argument applies at every order.
\end{proof}

\section{Closed forms, spectral calculus, and contraction semigroups}
\label{ap:formas-cerradas-calculo-espectral-semigrupos}

The elliptic operators appearing in global analysis are generally unbounded. Their domain is part of the definition and cannot be omitted. In $L^2(X)$, closed quadratic forms provide a particularly stable way to construct self-adjoint realizations; in $L^p(X)$, strongly continuous semigroups and their generators play the corresponding role. We collect here the results used in constructing the heat semigroup and Bessel potentials. For closed forms we follow \cite[Appendix~B]{Guneysu2017Covariant}, and for semigroups, \cite[Chapter~II]{EngelNagel2000}.

\subsection{Closed forms and the Friedrichs realization}

\begin{definition}\label{def:forma-sesquilineal-cerrada-semibounded}
Let $H$ be a complex Hilbert space. A \textbf{sesquilinear form} on $H$ is a map
\[
\mathfrak a\colon V\times V\longrightarrow\mathbb C,
\]
where $V\subseteq H$ is a dense vector subspace, linear in the first variable and antilinear in the second. The form is called:
\begin{enumerate}[label=(\alph*)]
\item \textbf{symmetric} if $\mathfrak a(u,v)=\overline{\mathfrak a(v,u)}$;
\item \textbf{nonnegative} if $\mathfrak a(u,u)\geq0$;
\item \textbf{bounded below} if there exists $c\in\mathbb R$ such that $\mathfrak a(u,u)\geq c\|u\|_H^2$;
\item \textbf{closed} if $V$ is complete with respect to the norm
\[
\|u\|_{\mathfrak a}^2
:=
\mathfrak a(u,u)+(1-c)\|u\|_H^2,
\]
where $c$ is any lower bound strictly less than $1$.
\end{enumerate}
The space $V$ is called the domain of the form and is denoted by $\mathcal D(\mathfrak a)$.
\end{definition}

Up to equivalence, the preceding norm does not depend on the chosen lower bound. For a nonnegative form we simply use
\[
\|u\|_{\mathfrak a}^2:=\|u\|_H^2+\mathfrak a(u,u).
\]

\begin{theorem}[First representation theorem for closed forms]\label{teo:representacion-formas-cerradas}\index{representation theorem for forms@representation theorem for closed forms}
Let $\mathfrak a$ be a densely defined, symmetric, nonnegative, closed sesquilinear form on a Hilbert space $H$. There exists a unique nonnegative self-adjoint operator
\[
A\colon \mathcal D(A)\subseteq H\longrightarrow H
\]
such that $\mathcal D(A)\subseteq\mathcal D(\mathfrak a)$ and
\[
\mathfrak a(u,v)=\langle Au,v\rangle_H,
\qquad
u\in\mathcal D(A),\quad v\in\mathcal D(\mathfrak a).
\]
More precisely,
\[
\mathcal D(A)
=
\left\{
 u\in\mathcal D(\mathfrak a)
 \middle|
 \text{there exists }f\in H\text{ such that }
 \mathfrak a(u,v)=\langle f,v\rangle_H
 \text{ for every }v\in\mathcal D(\mathfrak a)
\right\},
\]
and then $Au=f$.
\end{theorem}

\begin{proof}
Equip $V:=\mathcal D(\mathfrak a)$ with the inner product
\[
\langle u,v\rangle_{\mathfrak a}
:=
\langle u,v\rangle_H+\mathfrak a(u,v).
\]
Since the form is closed, $V$ is a Hilbert space. For each
$f\in H$, the antilinear functional
$v\mapsto\langle f,v\rangle_H$ is continuous on $V$, since
$\|v\|_H\leq\|v\|_{\mathfrak a}$. The antilinear version of
Theorem~\ref{teo:representacion-riesz-hilbert} provides a unique
$Rf\in V$ such that
\[
\langle Rf,v\rangle_H+\mathfrak a(Rf,v)
=
\langle f,v\rangle_H,
\qquad v\in V.
\]
The operator $R\colon H\longrightarrow H$, obtained by composing with the
inclusion $V\hookrightarrow H$, is linear and bounded. It is injective: if
$Rf=0$, then $\langle f,v\rangle_H=0$ for every $v\in V$, and
density of $V$ gives $f=0$. For $f,h\in H$, apply the defining
identity to $(f,Rh)$ and $(h,Rf)$; symmetry of $\mathfrak a$ gives
\[
 \langle f,Rh\rangle_H=\langle Rf,h\rangle_H.
\]
Thus $R$ is self-adjoint. Taking $v=Rf$ gives
\[
 \langle f,Rf\rangle_H
 =\|Rf\|_H^2+\mathfrak a(Rf,Rf)\geq0,
\]
so $R$ is positive. Its range is dense: if
$h\perp\operatorname{Ran}(R)$, self-adjointness implies $Rh=0$ and
injectivity gives $h=0$.

Define
\[
A:=R^{-1}-I,
\qquad
\mathcal D(A):=\operatorname{Ran}(R).
\]
If $u=Rf$, the identity defining $R$ becomes
\[
\mathfrak a(u,v)=\langle f-u,v\rangle_H
=
\langle Au,v\rangle_H,
\qquad v\in V.
\]
This formula proves the characterization of the domain and nonnegativity.
Moreover, $A$ is symmetric and $A+I=R^{-1}$ is surjective. If
$y\in\mathcal D(A^*)$, choose $x\in\mathcal D(A)$ with
\[
 (A+I)x=(A^*+I)y.
\]
Since $A\subseteq A^*$, we have $(A^*+I)(y-x)=0$. But
\[
 \ker(A^*+I)=\operatorname{Ran}(A+I)^\perp=\{0\},
\]
so $y=x\in\mathcal D(A)$. Thus $A=A^*$, without using the spectral
calculus introduced in the next subsection. Uniqueness
follows from the characterization of the domain.

\end{proof}

\begin{corollary}[Variational characterization of the resolvent]\label{cor:caracterizacion-variacional-resolvente-forma}\index{resolvent!variational characterization}
Let $A$ be the operator associated with the form $\mathfrak a$ in Theorem~\ref{teo:representacion-formas-cerradas}. If $\lambda>0$ and $f\in H$, then $u=(\lambda I+A)^{-1}f$ is the unique element of $\mathcal D(\mathfrak a)$ satisfying
\[
\mathfrak a(u,v)+\lambda\langle u,v\rangle_H
=
\langle f,v\rangle_H
\]
for every $v\in\mathcal D(\mathfrak a)$.
\end{corollary}

\begin{proof}
The expression $\langle u,v\rangle_{\mathfrak a,\lambda}:=\mathfrak a(u,v)+\lambda\langle u,v\rangle_H$ defines an inner product on $\mathcal D(\mathfrak a)$ whose norm is equivalent to the form norm. For $f\in H$, the functional $v\mapsto\langle f,v\rangle_H$ is continuous with respect to that norm. The Riesz representation theorem, Theorem~\ref{teo:representacion-riesz-hilbert}, provides a unique $u\in\mathcal D(\mathfrak a)$ satisfying the identity in the statement. Rewriting it as $\mathfrak a(u,v)=\langle f-\lambda u,v\rangle_H$ and using the characterization of $\mathcal D(A)$ in Theorem~\ref{teo:representacion-formas-cerradas} gives $u\in\mathcal D(A)$ and $Au=f-\lambda u$. Therefore $(\lambda I+A)u=f$, and the uniqueness already proved shows that $u=(\lambda I+A)^{-1}f$.
\end{proof}

\begin{definition}[Friedrichs realization]\label{def:realizacion-friedrichs}
Let $S$ be a symmetric, densely defined operator bounded
below on $H$, and choose $c\in\mathbb R$ such that
$S\geq cI$. The nonnegative form
\[
 \mathfrak b_c(u,v)
 :=\langle Su,v\rangle_H-c\langle u,v\rangle_H,
 \qquad u,v\in\mathcal D(S),
\]
is closable. Let $B_c$ be the operator associated with its closure by
Theorem~\ref{teo:representacion-formas-cerradas}. The operator
\[
 S_F:=B_c+cI,\qquad \mathcal D(S_F):=\mathcal D(B_c),
\]
is called the \textbf{Friedrichs realization} of $S$. It is independent of the
chosen lower bound $c$ and is a self-adjoint extension of $S$ with the
same lower bound.
\end{definition}

\begin{proof}
Complete $\mathcal D(S)$ with respect to the inner product
\[
(u,v)_c:=\langle u,v\rangle_H+\mathfrak b_c(u,v)
\]
and denote the resulting Hilbert space by $V_c$. The inclusion
into $H$ extends to a continuous operator $J_c\colon V_c\to H$.
We prove that it is injective. If $w\in\ker J_c$, choose
$u_j\in\mathcal D(S)$ such that $u_j\to w$ in $V_c$; then
$u_j\to0$ in $H$. For each $v\in\mathcal D(S)$, symmetry of $S$ gives
\[
(w,v)_c=\lim_{j\to\infty}(u_j,v)_c
=\lim_{j\to\infty}\langle u_j,(S+(1-c)I)v\rangle_H=0.
\]
Density of $\mathcal D(S)$ in $V_c$ implies $w=0$.
We may therefore identify $V_c$ with a subspace of $H$.
The limit of $\mathfrak b_c(u_j,v_j)$ for sequences converging
in $V_c$ defines a closed form on this subspace. It is precisely
the closure of $\mathfrak b_c$, since $\mathcal D(S)$ is dense in
the form norm.

If $u\in\mathcal D(S)$, the identity
$\mathfrak b_c(u,v)=\langle(S-cI)u,v\rangle_H$, initially for
$v\in\mathcal D(S)$, passes to the limit for every $v\in V_c$.
The characterization of the domain of the associated operator shows that
$B_cu=(S-cI)u$; thus $S_F$ extends $S$.
For two lower bounds $c,d$, the form norms are equivalent,
since their squares differ by $(d-c)\|u\|_H^2$ and both dominate
$\|u\|_H^2$. Their completions are therefore identified with the
same subspace of $H$. The closed forms satisfy
\[
\overline{\mathfrak b_c}(u,v)+c\langle u,v\rangle_H
=\overline{\mathfrak b_d}(u,v)+d\langle u,v\rangle_H.
\]
The variational characterization of their operators gives
$B_c+cI=B_d+dI$, including their domains. This proves
independence of the choice of lower bound.
\end{proof}

\begin{proposition}[Criterion for essential self-adjointness]\label{prop:criterio-autoadjuncion-esencial-semibounded}
Let $S$ be a symmetric, densely defined operator on a complex Hilbert space. Suppose that $S\geq0$. Then $S$ is essentially self-adjoint if and only if
\[
\ker(S^*+I)=\{0\}.
\]
More generally, if $S\geq-cI$, one may replace $I$ by $\lambda I$ with $\lambda>c$.
\end{proposition}

\begin{proof}
For $u\in\mathcal D(S)$,
\[
\|(S+I)u\|_H\|u\|_H
\geq
\operatorname{Re}\langle(S+I)u,u\rangle_H
\geq
\|u\|_H^2,
\]
so $\|(S+I)u\|_H\geq\|u\|_H$. The same estimate holds for the closure $\overline S$, so $\operatorname{Ran}(\overline S+I)$ is closed. On the other hand,
\[
\operatorname{Ran}(S+I)^\perp=\ker(S^*+I).
\]
If this kernel is trivial, $\operatorname{Ran}(S+I)$ is dense. The
preceding estimate also shows that
\[
 \operatorname{Ran}(\overline S+I)
 =\overline{\operatorname{Ran}(S+I)}=H.
\]
Let $y\in\mathcal D(S^*)=\mathcal D((\overline S)^*)$. Choose
$x\in\mathcal D(\overline S)$ such that
\[
 (\overline S+I)x=(S^*+I)y.
\]
Since $\overline S\subseteq S^*$, we obtain
$(S^*+I)(y-x)=0$. The hypothesis gives $y=x$, so
$\mathcal D(S^*)\subseteq\mathcal D(\overline S)$ and $\overline S$ is
self-adjoint. Conversely, if $\overline S$ is self-adjoint and
nonnegative, then $-1$ belongs to its resolvent set; hence
$\ker(S^*+I)=\ker(\overline S+I)=\{0\}$. For $S\geq-cI$, apply the
same argument to $S+cI$ and any shift $\lambda>c$.
\end{proof}

\subsection{Spectral calculus for self-adjoint operators}

\begin{definition}[Spectral measure]
\label{def:medida-espectral-proyecciones}
A \textit{spectral measure} on a measurable space $X$ is a map $E$
from its measurable sets to the orthogonal projections on $H$, such
that $E(X)=I$, $E(B\cap C)=E(B)E(C)$ and, for disjoint sets $B_j$,
\[
 E\Bigl(\bigcup_jB_j\Bigr)u=\sum_jE(B_j)u,\qquad u\in H,
\]
with convergence in norm. The scalar measure
$\mu_u(B):=\langle E(B)u,u\rangle_H=\|E(B)u\|_H^2$ is positive and
has total mass $\|u\|_H^2$.
\end{definition}

\begin{lemma}[Construction of the bounded spectral calculus]
\label{lem:calculo-espectral-autoadjunto-acotado}
Let $T\in\mathcal L(H)$ be self-adjoint and let $r>0$ satisfy
$r\geq\|T\|$. There is a unique spectral measure $E_T$ on
$I=[-r,r]$ representing $T$. For every Borel function
$m:I\to\mathbb C$, its integral $m(T)$ is a closed, densely defined
operator, with domain and norm given by
\[
 \mathcal D(m(T))=\{u\mid\int_I|m|^2\,d\mu_u<\infty\},\qquad
 \|m(T)u\|_H^2=\int_I|m|^2\,d\mu_u.
\]
Moreover, $m(T)^*=\overline m(T)$. For bounded functions the calculus
is an adjoint-preserving homomorphism and satisfies
$\|m(T)\|\leq\|m\|_\infty$.
\end{lemma}

\begin{proof}
We first construct the calculus for continuous functions. The
polarization identity used in the proof of
Theorem~\ref{teo:espectral-compacto-autoadjunto} does not require
compactness and gives
\[
 \|S\|=\sup_{\|u\|_H=1}|\langle Su,u\rangle_H|
 \quad\text{for every }S=S^*\in\mathcal L(H).
\]
We also need positivity of products of commuting positive operators.
If $P\geq0$ and $P\neq0$, set $C=I-\displaystyle\frac{P}{\|P\|}$.
Then $0\leq C\leq I$ and $\|C\|\leq1$. The binomial series
\[
 (1-z)^{\frac{1}{2}}=1-\sum_{j\geq1}a_jz^j,\qquad
 a_j=\frac{(2j)!}{(2j-1)(j!)^2\,4^j}>0,
 \qquad\sum_{j\geq1}a_j=1,
\]
defines $P^{\displaystyle\frac{1}{2}}=\|P\|^{\displaystyle\frac{1}{2}}(I-\sum_{j\geq1}a_jC^j)$ by absolute
convergence in operator norm. Powers of $C$ are positive: for even
exponents use $\langle C^{2k}u,u\rangle=\|C^ku\|^2$, and for odd
exponents apply positivity of $C$ to $C^ku$. Since $0\leq C^j\leq I$,
the partial sums defining the square root are positive. The binomial
identity, first for $|z|<1$ and then by absolute convergence, gives
$(P^{\displaystyle\frac{1}{2}})^2=P$. This square root commutes with every operator
commuting with $P$. For $P=0$ take the zero square root.
If $P,Q\geq0$ commute, then
\[
 PQ=P^{\frac{1}{2}}QP^{\frac{1}{2}}\geq0.
\]

Let $q$ be a real polynomial nonnegative on $I$. Unless it is zero,
its real factorization consists, after orienting the linear factors
and extracting a positive constant, of factors of the forms
\[
 (t-a)^{2k},\qquad t-a\ (a\leq-r),\qquad
 a-t\ (a\geq r),\qquad (t-a)^2+b^2\ (b\neq0).
\]
Interior roots have even order; endpoint roots are included among the
oriented linear factors. Each factor evaluated at $T$ is positive,
and all the factors commute. Hence $q(T)\geq0$.
For a complex polynomial $p$, set $M=\max_I|p|$ and apply this result
to $M^2-|p|^2$. We obtain
\[
 p(T)^*p(T)\leq M^2I,\qquad \|p(T)\|\leq\max_I|p|.
\]

Polynomials approximate every continuous function on $I$ uniformly.
To justify this step, rescale the interval to $[0,1]$. The Bernstein
polynomials of a continuous function $f$ satisfy
\[
 B_Nf(x)=\sum_{k=0}^N f(\frac{k}{N})\binom Nk x^k(1-x)^{N-k},
\]
\[
 |B_Nf(x)-f(x)|
 \leq\omega_f(\delta)+\frac{\|f\|_\infty}{2N\delta^2}.
\]
Here $\omega_f(\delta)$ is the modulus of continuity. Split the sum
into $|\displaystyle\frac{k}{N}-x|<\delta$ and its complement, and use that the mean and
variance of the binomial weights are $x$ and $\displaystyle\frac{x(1-x)}{N}$.
Choosing $\delta$ first and then $N$ proves uniform approximation;
real and imaginary parts are handled separately.
If $p_j\to f$ uniformly, the operator bound above makes $p_j(T)$
converge in norm to an operator $f(T)$ independent of the
approximation. This gives a unital, adjoint-preserving homomorphism,
which is positive because $f\geq0$ implies $f=(\sqrt f)^2$.

For $u\in H$, the positive functional
$f\mapsto\langle f(T)u,u\rangle_H$ has a unique positive regular
measure $\mu_u$ on $I$, by
Theorem~\ref{teo: representacion de Riesz en Cc}.
Its mass is $\|u\|_H^2$. Polarizing these measures gives complex
measures $\mu_{u,v}$ such that
\[
 \langle f(T)u,v\rangle_H=\int_I f\,d\mu_{u,v}.
\]
Uniqueness of the representation shows that $\mu_{u,v}(B)$ is
sesquilinear in $(u,v)$ and
$\mu_{u,u}(B)=\mu_u(B)\geq0$. Cauchy--Schwarz for this form gives
$|\mu_{u,v}(B)|^2\leq\mu_u(B)\mu_v(B)$.
Applying it on a finite partition and summing yields
\[
 \left|\int_I g\,d\mu_{u,v}\right|
 \leq\left(\int_I|g|^2\,d\mu_u\right)^{\frac{1}{2}}\|v\|_H
\]
for simple functions; uniform approximation extends it to bounded
Borel functions. In particular, the form $\mu_{u,v}(B)$ represents
a positive operator $E(B)$ with $0\leq E(B)\leq I$ and
$\langle E(B)u,v\rangle_H=\mu_{u,v}(B)$.

We show that $E(B)$ is a projection. Regularity of $\mu_u$ provides
$f_j\in C(I)$, $0\leq f_j\leq1$, converging to $\mathbf1_B$ in
$L^2(I,\mu_u)$. Indeed, choose compact $K_j$ and open $O_j$ with
$K_j\subset B\subset O_j$ and
$\mu_u(O_j\setminus K_j)\to0$, and continuous functions equal to one
on $K_j$ and zero outside $O_j$. The identity
\[
 \|(f_j-f_k)(T)u\|_H^2=\int_I|f_j-f_k|^2\,d\mu_u
\]
shows that $f_j(T)u$ converges in norm. The estimate for $\mu_{u,v}$
identifies its limit with $E(B)u$. Consequently,
\[
 \|E(B)u\|_H^2=\mu_u(B)=\langle E(B)u,u\rangle_H,
\]
and polarization gives $E(B)^2=E(B)$. If $B,C$ are disjoint, scalar
additivity gives $E(B\cup C)=E(B)+E(C)$.
Since this sum is a projection, its two terms are orthogonal:
evaluating $E(B)E(C)+E(C)E(B)=0$ on the range of $E(B)$ yields
$\|E(C)E(B)u\|_H^2=0$. Decomposing two sets into their intersection
and differences proves $E(B)E(C)=E(B\cap C)$.
Countable additivity is strong because the squared norm of the
remainder applied to $u$ is the $\mu_u$ mass of the remaining union.
Thus $E$ is a spectral measure.

For $m=\sum_{j=1}^Nc_j\mathbf1_{B_j}$, with disjoint $B_j$, define
$m(T)=\sum_jc_jE(B_j)$. Orthogonality gives the norm formula,
the product rules, and $m(T)^*=\overline m(T)$.
Uniform approximation by simple functions defines the calculus for
all bounded Borel functions. For continuous functions it agrees with
the previous construction because both representations have the same
inner products. In particular, $T=\int_I t\,dE(t)$.
Uniqueness follows because any representing spectral measure has the
same polynomial integrals and, by approximation, the same continuous
function integrals. Riesz uniqueness identifies its scalar measures,
and polarization identifies its projections.

Finally, for a finite Borel function $m$, set
$F_N=E(\{|m|\leq N\})$. The bounded operators
$(m\mathbf1_{\{|m|\leq N\}})(T)$ applied to $u$ converge in norm
exactly when $\int_I|m|^2\,d\mu_u<\infty$. Their limit defines
$m(T)u$ and satisfies the norm formula. The domain contains every
range of $F_N$ and is dense because $F_N\to I$ strongly.
If $u_j\to u$ and $m(T)u_j\to v$, then
\[
 F_Nv=(m\mathbf1_{\{|m|\leq N\}})(T)u.
\]
These norms are bounded by $\|v\|_H$. Monotone convergence gives
$u\in\mathcal D(m(T))$, and then $m(T)u=v$. Thus the operator is
closed. Truncation proves $\overline m(T)\subseteq m(T)^*$.
If $v\in\mathcal D(m(T)^*)$ and $h=m(T)^*v$, testing the adjoint
identity on $F_NH$ gives
\[
 (\overline m\mathbf1_{\{|m|\leq N\}})(T)v=F_Nh.
\]
Monotone convergence again implies
$v\in\mathcal D(\overline m(T))$ and $\overline m(T)v=h$.
This proves equality of the adjoints. To complete the multiplication
rules, projections commute with every truncation; thus, for
$u\in\mathcal D(m(T))$,
\[
 \mu_{m(T)u}(B)=\|E(B)m(T)u\|_H^2
               =\int_B|m|^2\,d\mu_u.
\]
For bounded $m$ this follows directly from simple functions; for
unbounded $m$ take the norm limit of the truncations. This formula
determines when $m(T)u$ belongs to another multiplier's domain and
justifies the products used below.
\end{proof}

\begin{theorem}[Spectral theorem]\label{teo:calculo-espectral-autoadjunto-no-negativo}\index{spectral theorem!unbounded self-adjoint operators}
Let $A$ be a self-adjoint operator on a complex Hilbert space $H$. There exists a unique spectral measure $E_A$ on $\mathbb R$ such that
\[
A=\int_{\mathbb R}\lambda\,dE_A(\lambda).
\]
If $m\colon \mathbb R\longrightarrow\mathbb C$ is Borel measurable, define
\[
m(A):=\int_{\mathbb R}m(\lambda)\,dE_A(\lambda)
\]
with domain
\[
\mathcal D(m(A))
=
\left\{
 u\in H
 \middle|
 \int_{\mathbb R}|m(\lambda)|^2\,d\mu_u(\lambda)<\infty
\right\},
\]
where $\mu_u(B):=\langle E_A(B)u,u\rangle_H$. We have
\[
\|m(A)u\|_H^2
=
\int_{\mathbb R}|m(\lambda)|^2\,d\mu_u(\lambda).
\]
If $m$ is bounded, then $m(A)\in\mathcal L(H)$ and
\[
\|m(A)\|_{\mathcal L(H)}
\leq\|m\|_{\infty,\mathbb R}
:=\sup_{\lambda\in\mathbb R}|m(\lambda)|.
\]
Moreover, the calculus respects sums, products, and conjugation on their natural domains.
\end{theorem}

\begin{proof}
We reduce the unbounded operator to a bounded self-adjoint operator,
without using spectral calculus for $A$ at this stage.
Self-adjointness and the identity
\[
 \|(A\pm iI)u\|_H^2=\|Au\|_H^2+\|u\|_H^2
 \qquad(u\in\mathcal D(A))
\]
show that $A\pm iI$ have closed ranges and zero kernels.
Their ranges are dense because their orthogonal complements are
$\ker(A\mp iI)=\{0\}$. Their inverses
$Q_\pm=(A\pm iI)^{-1}$ therefore exist and are bounded.
The resolvent identity, obtained by subtracting the equations for
the two inverses, gives $Q_+Q_-=Q_-Q_+$; moreover,
$Q_+^*=Q_-$. Define
\[
 R:=Q_-Q_+=(I+A^2)^{-1}.
\]
The last equality includes the domains: for $u=Q_-Q_+f$ we have
$Au=Q_+f+iu\in\mathcal D(A)$ and $(I+A^2)u=f$; the inverse
identity is verified on $\mathcal D(A^2)$.
Thus $R$ is positive, injective, $\|R\|\leq1$, and its range is
$\mathcal D(A^2)$.

Equip $V=\mathcal D(A)$ with its graph norm
$\|u\|_V^2=\|u\|_H^2+\|Au\|_H^2$. Since $A$ is closed, $V$ is
complete. For $f\in H$ and $v\in V$, self-adjointness gives
\[
 \langle Rf,v\rangle_H+\langle ARf,Av\rangle_H
 =\langle f,v\rangle_H.
\]
Consequently, $\operatorname{Ran}R$ is dense in $V$: a vector
orthogonal to that range for the inner product of $V$ would, by
this identity, be orthogonal to all of $H$.

Apply Lemma~\ref{lem:calculo-espectral-autoadjunto-acotado} to $R$.
Its measure is concentrated on $[0,1]$, since a spectral vector on
a negative interval would contradict $R\geq0$.
Also $E_R(\{0\})=0$ by injectivity.
Define
\[
 W:=\mathcal D(R^{-\frac{1}{2}}),\qquad
 \|u\|_W:=\|R^{-\frac{1}{2}}u\|_H.
\]
Use the function $t^{-\displaystyle\frac{1}{2}}$ on $(0,1]$, with any value at zero.
The bounded calculus and its extension by truncation give
$W=\operatorname{Ran}R^{\displaystyle\frac{1}{2}}$: for $u\in W$ use
$u=R^{\displaystyle\frac{1}{2}}(R^{-\displaystyle\frac{1}{2}}u)$, and verify the inverse identity in the same
way. Hence $W$ is complete and $\|u\|_H\leq\|u\|_W$.
The range of $R$ is dense in $W$ as well. For $u\in W$, the vectors
$E_R([\displaystyle\frac{1}{N},1])u$ belong to $\operatorname{Ran}R$ and converge to $u$
in this norm because the integral of $t^{-1}$ over $(0,\displaystyle\frac{1}{N})$ tends
to zero. On that range the two norms agree:
\[
 \|Rf\|_W^2=\langle Rf,f\rangle_H=\|Rf\|_V^2.
\]
Their completions are identified inside $H$. Therefore,
\[
 \operatorname{Ran}R^{\frac{1}{2}}=\mathcal D(A),\qquad
 \|AR^{\frac{1}{2}}f\|_H^2=\|f\|_H^2-\|R^{\frac{1}{2}}f\|_H^2.
\]

Set $S:=AR^{\displaystyle\frac{1}{2}}\in\mathcal L(H)$. The last identity yields
$\|S\|\leq1$. This operator is self-adjoint. Indeed, $R$ preserves
$\mathcal D(A)$ and commutes there with $A$, as follows from
$AQ_\pm u=Q_\pm Au=u\mp iQ_\pm u$.
Uniform polynomial approximation of the square root and closedness
of $A$ show that $R^{\displaystyle\frac{1}{2}}$ preserves this domain and
$AR^{\displaystyle\frac{1}{2}}u=R^{\displaystyle\frac{1}{2}}Au$ for $u\in\mathcal D(A)$.
For $f,g$ in that domain,
\[
 \langle Sf,g\rangle_H
 =\langle R^{\frac{1}{2}}f,Ag\rangle_H
 =\langle f,Sg\rangle_H.
\]
Density and continuity extend the identity to $H$, so $S=S^*$.
Polarizing the norm identity gives
\[
 S^2=I-R.
\]

Let $E_S$ be the measure of the bounded operator $S$. It is
concentrated on $[-1,1]$, and $E_S(\{-1,1\})=0$ because
$I-S^2=R$ is injective.
Compatibility of continuous calculus with composition, obtained by
polynomial approximation, gives $R^{\displaystyle\frac{1}{2}}=(1-S^2)^{\displaystyle\frac{1}{2}}$.
On $(-1,1)$ set
\[
 b(s):=\frac{s}{\sqrt{1-s^2}},\qquad
 h(t):=\frac{t}{\sqrt{1+t^2}}.
\]
These are inverse homeomorphisms between $(-1,1)$ and $\mathbb R$.
The endpoint values of $b$ do not affect its operator.
Since $1+|b(s)|^2=(1-s^2)^{-1}$, the lemma gives
\[
 \mathcal D(b(S))=\operatorname{Ran}(1-S^2)^{\frac{1}{2}}
 =\mathcal D(A),\qquad
 b(S)R^{\frac{1}{2}}f=Sf=AR^{\frac{1}{2}}f.
\]
Thus $A=b(S)$. For a Borel set $B\subset\mathbb R$, define
\[
 E_A(B):=E_S(b^{-1}(B)\cap(-1,1)).
\]
This is a spectral measure representing $A$.
For Borel $m$, define its integral by $(m\circ b)(S)$.
Changing variables in the scalar measures gives precisely the domain
and norm formula in the statement. The bound for bounded $m$ is
the bound in the lemma.

We specify the rules for unbounded functions.
We have $m(A)^*=\overline m(A)$ and
\[
 (m+n)(A)u=m(A)u+n(A)u
 \quad\bigl(u\in\mathcal D(m(A))\cap\mathcal D(n(A))\bigr).
\]
The spectral sum may have a larger domain than this intersection.
Moreover, orthogonality of the projections and truncation give
$d\mu_{n(A)u}=|n|^2\,d\mu_u$ for $u\in\mathcal D(n(A))$.
Consequently,
\[
 \begin{split}
 \mathcal D(m(A)n(A))
 &=\mathcal D(n(A))\cap\mathcal D((mn)(A)),\\
 m(A)n(A)u&=(mn)(A)u
 \quad\text{on this domain}.
 \end{split}
\]
First verify these identities on sets where both multipliers are
bounded. The norm formula then allows passage to the limit on
each indicated domain.

For uniqueness, any measure representing $A$ yields, through its
integral calculus,
\[
 (I+A^2)^{-1}=\int_{\mathbb R}\frac{1}{1+t^2}\,dE(t),\qquad
 S=\int_{\mathbb R}h(t)\,dE(t).
\]
Its image under $h$ represents $S$ on $[-1,1]$.
Uniqueness in the lemma identifies that image with $E_S$.
Since $h$ is a homeomorphism onto $(-1,1)$ and there is no mass at
the endpoints, we conclude that $E=E_A$.
\end{proof}

Sources: \cite{Guneysu2017Covariant,reed_simon_vol1}.

\begin{corollary}[Square-root domain associated with a closed form]
\label{cor:dominio-raiz-forma-cerrada}
\index{closed form!square-root domain}
Let $\mathfrak a$ be a form as in
Theorem~\ref{teo:representacion-formas-cerradas}, and let $A$ be its associated self-adjoint
operator. Then
\[
 \mathcal D(A^{\frac{1}{2}})=\mathcal D(\mathfrak a),
 \qquad
 \mathfrak a(u,v)
 =\langle A^{\frac{1}{2}}u,A^{\frac{1}{2}}v\rangle_H.
\]
\end{corollary}
\begin{proof}
Return to $V=\mathcal D(\mathfrak a)$ and the resolvent $R$ constructed in
the proof of Theorem~\ref{teo:representacion-formas-cerradas}. The
subspace $\mathcal D(A)=\operatorname{Ran}(R)$ is dense in $V$ for the
form norm. Indeed, if $w\in V$ is orthogonal to
$\operatorname{Ran}(R)$ with respect to
$\langle\cdot,\cdot\rangle_{\mathfrak a}$, then, for every $f\in H$,
\[
 0=\langle Rf,w\rangle_{\mathfrak a}=\langle f,w\rangle_H,
\]
and hence $w=0$.

For $u\in\mathcal D(A)$, the spectral calculus and the representation of the
form give
\[
 \mathfrak a(u,u)=\langle Au,u\rangle_H
 =\|A^{\frac{1}{2}}u\|_H^2.
\]
Moreover, $\mathcal D(A)$ is a core for $A^{\displaystyle\frac{1}{2}}$. If
$u\in\mathcal D(A^{\displaystyle\frac{1}{2}})$ and
$u_N:=E_A([0,N])u$, then $u_N\in\mathcal D(A)$ and
\[
 \|u_N-u\|_H^2
 +\|A^{\frac{1}{2}}(u_N-u)\|_H^2
 =
 \int_{(N,+\infty)}(1+\lambda)\,d\mu_u(\lambda)
 \longrightarrow0.
\]
Thus $V$ and $\mathcal D(A^{\displaystyle\frac{1}{2}})$ are the completions of
$\mathcal D(A)$ with respect to the same norm. They agree as subspaces of
$H$, and the sesquilinear identity follows by polarization.
\end{proof}

\begin{corollary}\label{cor:potencias-espectrales-operador-no-negativo}
Let $A$ be self-adjoint and nonnegative. For $s\in\mathbb R$, define $(I+A)^s$ and, for $s\geq0$, $A^s$ by the spectral calculus. If $s>0$, then $(I+A)^{-s}$ is bounded, injective, and
\[
\operatorname{Ran}\bigl((I+A)^{-s}\bigr)=\mathcal D\bigl((I+A)^s\bigr).
\]
The operators $(I+A)^s$ and $(I+A)^{-s}$ are inverses of one another on these domains.
\end{corollary}

\begin{proof}
Let $E_A$ be the spectral measure of $A$ and let
$\mu_u(B):=\|E_A(B)u\|_H^2$. Since $A$ is nonnegative, these measures
are concentrated on $[0,\infty)$. The multipliers
$(1+\lambda)^s$, with $s\in\mathbb R$, and $\lambda^s$, with $s\geq0$,
are Borel and finite on this interval; for $s=0$ we take
$\lambda^0=1$, including at $\lambda=0$. By
Theorem~\ref{teo:calculo-espectral-autoadjunto-no-negativo}, they define
the operators in the statement on their spectral domains.

Fix $s>0$ and put $B:=(I+A)^s$ and $C:=(I+A)^{-s}$.
The inequality $(1+\lambda)^{-s}\leq1$ implies that $C$ is bounded
and $\|C\|\leq1$. Moreover,
\[
 \|Cu\|_H^2
 =\int_{[0,\infty)}(1+\lambda)^{-2s}\,d\mu_u(\lambda).
\]
The integrand is strictly positive. If $Cu=0$, its zero integral
forces $\mu_u$ to be the zero measure; its mass is
$\|u\|_H^2$, so $u=0$. Thus $C$ is injective.

The rule for the scalar measures of a transformed vector gives
$d\mu_{Cu}=(1+\lambda)^{-2s}\,d\mu_u$. Consequently,
\[
 \int_{[0,\infty)}(1+\lambda)^{2s}\,d\mu_{Cu}(\lambda)
 =\|u\|_H^2<\infty.
\]
This proves that $Cu\in\mathcal D(B)$ for every $u\in H$, and the
product rule on this domain gives $BCu=u$. Conversely,
if $v\in\mathcal D(B)$, boundedness of $C$ allows us to apply
the same rule and obtain $CBv=v$. In particular,
$v=C(Bv)\in\operatorname{Ran}(C)$. The two inclusions prove
$\operatorname{Ran}(C)=\mathcal D(B)$ and both inverse identities
with their domains.
\end{proof}

\begin{proposition}[Spectral semigroup]\label{prop:semigrupo-espectral-autoadjunto}
Let $A$ be self-adjoint and nonnegative. Then $T_t:=e^{-tA}$, $t\geq0$, is a strongly continuous semigroup of self-adjoint contractions. For each $u\in H$, $T_tu\in\mathcal D(A^m)$ if $t>0$ and $m\in\mathbb N$, and
\[
\frac{d}{dt}T_tu=-AT_tu
\qquad(t>0).
\]
If $u\in\mathcal D(A)$, the right derivative at $t=0$ also exists and equals $-Au$.
\end{proposition}

\begin{proof}
The algebraic properties follow from $e^{-(t+s)\lambda}=e^{-t\lambda}e^{-s\lambda}$, and the bound $|e^{-t\lambda}|\leq1$ gives contractivity. For strong continuity, observe that
\[
\|T_tu-u\|_H^2
=
\int_{[0,+\infty)}|e^{-t\lambda}-1|^2\,d\mu_u(\lambda).
\]
The integrand converges pointwise to zero as $t\to0^+$ and is dominated by the constant $4$, which is integrable with respect to the finite measure $\mu_u$; dominated convergence proves $T_tu\to u$.

If $t>0$, then
\[
\|A^mT_tu\|_H^2
=
\int_{[0,+\infty)}\lambda^{2m}e^{-2t\lambda}\,d\mu_u(\lambda)
\leq
\left(\sup_{\lambda\geq0}\lambda^{2m}e^{-2t\lambda}\right)\|u\|_H^2,
\]
so $T_tu\in\mathcal D(A^m)$. Now fix $t_0>0$ and take $|h|<\displaystyle\frac{t_0}{2}$. The mean value theorem gives
\[
\left|
\frac{e^{-(t_0+h)\lambda}-e^{-t_0\lambda}}{h}
\right|
\leq
\lambda e^{-\frac{t_0}{2}\lambda}.
\]
Thus the difference between this quotient and $-\lambda e^{-t_0\lambda}$ is dominated by $2\lambda e^{-\displaystyle\frac{t_0}{2}\lambda}$, a bounded function of $\lambda$. Dominated convergence with respect to $\mu_u$ proves differentiability at $t_0$ and identity $\displaystyle\frac{d}{dt}T_tu=-AT_tu$.

Finally, if $u\in\mathcal D(A)$, then
\[
\left|
\frac{e^{-t\lambda}-1}{t}
\right|
\leq\lambda
\]
for $t>0$. The square of the dominating function is integrable with respect to $\mu_u$, precisely because $u\in\mathcal D(A)$. One last application of dominated convergence gives the right derivative at $t=0$.
\end{proof}

\subsection{Strongly continuous semigroups on Banach spaces}

\begin{definition}\label{def:semigrupo-C0-contraccion}\index{strongly continuous semigroup}
Let $X$ be a Banach space. A family $(T_t)_{t\geq0}\subseteq\mathcal L(X)$ is a \textbf{strongly continuous semigroup}, or a $C_0$-semigroup, if
\[
T_0=I,
\qquad
T_{t+s}=T_tT_s,
\qquad
\lim_{t\to0^+}\|T_tx-x\|_X=0
\]
for every $x\in X$. It is a \textbf{contraction semigroup} if $\|T_t\|\leq1$ for every $t\geq0$.
\end{definition}

\begin{definition}\label{def:generador-semigrupo-C0-apendice}
The \textbf{infinitesimal generator} of a $C_0$-semigroup $(T_t)$ is the operator
\[
Gx:=\lim_{t\to0^+}\frac{T_tx-x}{t}
\]
defined on the set $\mathcal D(G)$ of elements for which the limit exists in $X$.
\end{definition}

\begin{proposition}\label{prop:propiedades-generador-semigrupo-C0}
Let $G$ be the generator of a $C_0$-semigroup $(T_t)_{t\geq0}$ on $X$. Then $G$ is closed and densely defined, and $\mathcal D(G^m)$ is dense in $X$ for every $m\in\mathbb N$. Moreover, $T_t\mathcal D(G)\subseteq\mathcal D(G)$, $GT_tx=T_tGx$ for $x\in\mathcal D(G)$, and
\[
T_tx-x=\int_0^tT_sGx\,ds,
\qquad x\in\mathcal D(G).
\]
\end{proposition}

\begin{proof}
We first prove local boundedness in operator norm. If there were no
$\delta>0$ and $M<\infty$ with $\sup_{0\leq t\leq\delta}\|T_t\|\leq M$,
we could choose $0<t_j<\displaystyle\frac{1}{j}$ with $\|T_{t_j}\|>j$. For each $x\in X$,
the sequence $T_{t_j}x$ converges to $x$ and is therefore bounded.
Theorem~\ref{teo: banach steinhaus} would give
$\sup_j\|T_{t_j}\|<\infty$, a contradiction. Take $M\geq1$.
Applying the semigroup property to $t=k\delta+r$ with
$0\leq r<\delta$ gives, for each $R>0$, a constant $C_R$ such that
\[
 \sup_{0\leq t\leq R}\|T_t\|\leq C_R.
\]
Strong right continuity at every $t\geq0$ follows from
$T_{t+h}x-T_tx=T_t(T_hx-x)$. For $0<h<t$, the identity
\[
 T_tx-T_{t-h}x=T_{t-h}(T_hx-x)
\]
and the local bound prove left continuity. Thus
$t\mapsto T_tx$ is continuous on $[0,\infty)$ for every $x\in X$.

If $x\in\mathcal D(G)$, then
\[
 \frac{T_hT_tx-T_tx}{h}
 =T_t\frac{T_hx-x}{h}\longrightarrow T_tGx.
\]
Hence $T_tx\in\mathcal D(G)$ and $GT_tx=T_tGx$.
The same quotient gives the right derivative of $T_tx$.
For $t>0$, the left derivative follows from
\[
 \frac{T_tx-T_{t-h}x}{h}
 =T_{t-h}\frac{T_hx-x}{h}\longrightarrow T_tGx;
\]
this limit uses local boundedness and strong continuity on $Gx$.
The derivative $T_tGx$ is continuous. The fundamental theorem of calculus gives
\[
 T_tx-x=\int_0^tT_sGx\,ds,
 \qquad x\in\mathcal D(G).
\]

For $h>0$, define the bounded operator
$S_hx:=h^{-1}\int_0^hT_sx\,ds$. Continuity at zero implies
$S_hx\to x$ as $h\to0^+$. For $t>0$,
\[
 \frac{T_tS_hx-S_hx}{t}
 =\frac{1}{h}\left(
 \frac{1}{t}\int_h^{h+t}T_sx\,ds
 -\frac{1}{t}\int_0^tT_sx\,ds\right)
 \longrightarrow\frac{T_hx-x}{h}.
\]
Thus $S_h(X)\subseteq\mathcal D(G)$ and
$GS_hx=\displaystyle\frac{T_hx-x}{h}$. Moreover, $S_hT_t=T_tS_h$, so
$GS_hx=S_hGx$ for $x\in\mathcal D(G)$. Induction gives
$S_h\mathcal D(G^k)\subseteq\mathcal D(G^{k+1})$ and
$S_h^m(X)\subseteq\mathcal D(G^m)$. If $0<h\leq\delta$, then
$\|S_h\|\leq M$ and
\[
 \|S_h^mx-x\|
 \leq\sum_{j=0}^{m-1}\|S_h\|^j\|S_hx-x\|
 \leq\left(\sum_{j=0}^{m-1}M^j\right)\|S_hx-x\|
 \longrightarrow0.
\]
This proves density of $\mathcal D(G^m)$, including $m=1$.

Finally, let $x_n\in\mathcal D(G)$ satisfy $x_n\to x$ and $Gx_n\to y$.
The integral identity already proved and the bound $C_t$ allow passage to the limit:
\[
 T_tx-x=\int_0^tT_sy\,ds,
 \qquad
 \left\|\int_0^tT_s(Gx_n-y)\,ds\right\|
 \leq tC_t\|Gx_n-y\|\longrightarrow0.
\]
Dividing by $t$ and letting $t\to0^+$, strong continuity on $y$
gives $x\in\mathcal D(G)$ and $Gx=y$. Therefore $G$ is closed.
\end{proof}

\begin{definition}\label{def:disipativo-acretivo-dualidad}
Let $A\colon \mathcal D(A)\subseteq X\longrightarrow X$ be densely defined. The operator $A$ is called \textbf{dissipative} if for each $x\in\mathcal D(A)$ there exists $x^*\in X'$ such that
\[
\|x^*\|_{X'}=\|x\|_X,
\qquad
x^*(x)=\|x\|_X^2,
\qquad
\operatorname{Re}x^*(Ax)\leq0.
\]
It is called \textbf{accretive} if $-A$ is dissipative. A dissipative operator is \textbf{maximally dissipative} if $\operatorname{Ran}(\lambda I-A)=X$ for some, and hence every, $\lambda>0$.
\end{definition}

\begin{proposition}[Resolvent as a Laplace transform]\label{prop:resolvente-transformada-laplace-semigrupo}
Let $(T_t)_{t\geq0}$ be a contraction $C_0$-semigroup with generator $G$. For each $\lambda>0$ and $x\in X$, the Bochner integral
\[
R_\lambda x:=\int_0^\infty e^{-\lambda t}T_tx\,dt
\]
converges, and
\[
R_\lambda=(\lambda I-G)^{-1},
\qquad
\|R_\lambda\|\leq\frac{1}{\lambda}.
\]
\end{proposition}

\begin{proof}
Strong continuity and the estimate
\[
 \int_0^\infty e^{-\lambda t}\|T_tx\|\,dt
 \leq\lambda^{-1}\|x\|
\]
prove existence of the Bochner integral and the bound on $R_\lambda$.
If $x\in\mathcal D(G)$, the preceding proposition gives
$\displaystyle\frac{d}{dt}T_tx=T_tGx$. Integration by parts on $[0,a]$ yields
\[
 \int_0^a e^{-\lambda t}T_t(\lambda I-G)x\,dt
 =x-e^{-\lambda a}T_ax.
\]
The last term tends to zero, since
$\|e^{-\lambda a}T_ax\|\leq e^{-\lambda a}\|x\|$.
Consequently, $R_\lambda(\lambda I-G)x=x$.

To prove the other identity, fix $x\in X$ and $a>0$ and set
$F_ax:=\int_0^a e^{-\lambda s}T_sx\,ds$. The semigroup property gives
\[
 T_hF_ax=e^{\lambda h}\int_h^{a+h}e^{-\lambda r}T_rx\,dr.
\]
Subtracting $F_ax$, dividing by $h$, and using continuity of the integrand
at $0$ and $a$, we obtain
\[
 \lim_{h\to0^+}\frac{T_hF_ax-F_ax}{h}
 =\lambda F_ax+e^{-\lambda a}T_ax-x.
\]
Hence $F_ax\in\mathcal D(G)$ and
$GF_ax=\lambda F_ax+e^{-\lambda a}T_ax-x$.
As $a\to\infty$, we have $F_ax\to R_\lambda x$ and
$GF_ax\to\lambda R_\lambda x-x$. Since $G$ is closed,
$R_\lambda x\in\mathcal D(G)$ and
$(\lambda I-G)R_\lambda x=x$. The two identities prove
$R_\lambda=(\lambda I-G)^{-1}$.
\end{proof}

\begin{theorem}[Lumer--Phillips]\label{teo:lumer-phillips}\index{Lumer--Phillips theorem}
Let $A$ be a closed, densely defined linear operator on a Banach space $X$. The following assertions are equivalent:
\begin{enumerate}[label=(\alph*)]
\item $A$ is the generator of a contraction $C_0$-semigroup;
\item $A$ is maximally dissipative;
\item $A$ is dissipative and $\operatorname{Ran}(\lambda I-A)=X$ for some $\lambda>0$.
\end{enumerate}
In this case,
\[
\|(\lambda I-A)^{-1}\|_{\mathcal L(X)}\leq\frac{1}{\lambda},
\qquad \lambda>0.
\]
\end{theorem}

\begin{proof}
First suppose that $A$ generates a contraction semigroup.
For $x\in\mathcal D(A)$, choose by Hahn--Banach an element $x^*\in X'$
with $\|x^*\|=\|x\|$ and $x^*(x)=\|x\|^2$. Since
$\operatorname{Re}x^*(T_tx)\leq\|x\|^2$, dividing the difference
by $t>0$ and passing to the limit gives $\operatorname{Re}x^*(Ax)\leq0$.
Proposition~\ref{prop:resolvente-transformada-laplace-semigrupo}
proves surjectivity of $\lambda I-A$ for every $\lambda>0$.

Now suppose (c). Dissipativity implies, for
$x\in\mathcal D(A)$ and $\lambda>0$,
\[
\lambda\|x\|^2
\leq\operatorname{Re}x^*((\lambda I-A)x)
\leq\|x\|\,\|(\lambda I-A)x\|.
\]
Thus $\lambda I-A$ is injective and its inverse on its range has
norm at most $\displaystyle\frac{1}{\lambda}$. The set $I$ of positive parameters
for which it is surjective is nonempty by hypothesis. To prove
openness, fix $\lambda\in I$ and set $R_\lambda:=(\lambda I-A)^{-1}$.
If $|\mu-\lambda|<\displaystyle\frac{\lambda}{2}$, then $\mu>0$ and
$\|(\mu-\lambda)R_\lambda\|<\displaystyle\frac{1}{2}$. On $\mathcal D(A)$,
\[
 \mu I-A=[I+(\mu-\lambda)R_\lambda](\lambda I-A).
\]
The first factor has a bounded inverse given by the convergent series
$\sum_{j=0}^\infty[-(\mu-\lambda)R_\lambda]^j$.
Hence $\mu I-A$ is surjective and
$(\displaystyle\frac{\lambda}{2},3\displaystyle\frac{\lambda}{2})\subseteq I$, proving openness.
It is also closed in $(0,\infty)$:
if $\lambda_j\in I$ and $\lambda_j\to\lambda>0$, the resolvent
identity and its bounds show that
$R_{\lambda_j}:=(\lambda_j I-A)^{-1}$ is Cauchy in norm,
since
\[
\|R_{\lambda_j}-R_{\lambda_k}\|
\leq\frac{|\lambda_j-\lambda_k|}{\lambda_j\lambda_k}.
\]
If its limit is $R$, for $y\in X$ we have
$R_{\lambda_j}y\to Ry$ and
$AR_{\lambda_j}y=\lambda_jR_{\lambda_j}y-y\to\lambda Ry-y$.
Closedness of $A$ gives $(\lambda I-A)Ry=y$.
Since $(0,\infty)$ is connected, $I=(0,\infty)$.

For each $\lambda>0$, consider the Yosida approximants
\[
A_\lambda:=\lambda A(\lambda I-A)^{-1}.
\]
They are bounded, since $A_\lambda=\lambda^2R_\lambda-\lambda I$, and
\[
\|e^{tA_\lambda}\|
\leq e^{-\lambda t}e^{t\lambda^2\|R_\lambda\|}\leq1,
\qquad t\geq0.
\]
Moreover, $\lambda R_\lambda y\to y$ for every $y\in X$.
Indeed, for $y\in\mathcal D(A)$,
\[
\lambda R_\lambda y-y=R_\lambda Ay,
\qquad \|R_\lambda Ay\|\leq\lambda^{-1}\|Ay\|,
\]
and density of $\mathcal D(A)$, together with
$\|\lambda R_\lambda\|\leq1$, extends the limit to every $X$.
Thus $A_\lambda x=\lambda R_\lambda Ax\to Ax$ if
$x\in\mathcal D(A)$.

The resolvent identity shows that $A_\lambda$ and $A_\mu$
commute. Differentiating
$s\mapsto e^{(t-s)A_\lambda}e^{sA_\mu}x$ and integrating gives
\[
\|e^{tA_\lambda}x-e^{tA_\mu}x\|
\leq t\|(A_\lambda-A_\mu)x\|,
\qquad x\in\mathcal D(A).
\]
Thus the contractions converge strongly, uniformly for
$0\leq t\leq T$, on the dense subspace $\mathcal D(A)$.
The uniform bound by one extends this convergence to all of $X$.
Denote the limit by $T_t$. Passing to the limit in
$e^{(t+s)A_\lambda}=e^{tA_\lambda}e^{sA_\lambda}$ gives the semigroup
property. For $x\in\mathcal D(A)$,
$\|T_tx-x\|\leq t\|Ax\|$; density and contractivity give
strong continuity at zero for every $x\in X$.

If $G$ is its generator, the identity
\[
e^{tA_\lambda}x-x=\int_0^t e^{sA_\lambda}A_\lambda x\,ds
\]
passes to the limit, for $x\in\mathcal D(A)$, as
$\displaystyle T_tx-x=\int_0^tT_sAx\,ds$. Dividing by $t$ gives
$Gx=Ax$, so $A\subseteq G$.
For $y\in\mathcal D(G)$, choose $x\in\mathcal D(A)$ with
$(\lambda I-A)x=(\lambda I-G)y$. Then
$(\lambda I-G)(y-x)=0$. The semigroup resolvent formula
implies $y=x$, and hence $G=A$. This proves (a).
The argument also proves that surjectivity for
one parameter is equivalent to surjectivity for all parameters, completing
the equivalence with (b) and the resolvent bound.
\end{proof}

\begin{proposition}[Exponential formula]\label{prop:formula-exponencial-semigrupo}
If $G$ generates a contraction $C_0$-semigroup $(T_t)$, then
\[
T_tx
=
\lim_{n\to\infty}
\left(I-\frac{t}{n}G\right)^{-n}x
\]
for every $x\in X$, uniformly for $t$ in compact intervals of $[0,\infty)$.
\end{proposition}

\begin{proof}
For $h>0$, the Laplace formula for the resolvent gives
\[
(I-hG)^{-1}x=\int_0^\infty e^{-r}T_{hr}x\,dr.
\]
Iterating this identity, Fubini and the semigroup property
yield
\[
(I-hG)^{-n}x
=\frac{1}{(n-1)!}\int_0^\infty e^{-r}r^{n-1}T_{hr}x\,dr.
\]
To verify the factor, the step from $n$ to $n+1$ integrates over
$0<s<r$ and uses $\displaystyle \int_0^r s^{n-1}\,ds=\displaystyle\frac{r^n}{n}$.
The measure $d\nu_n(r):=e^{-r}r^{n-1}\displaystyle\frac{dr}{(n-1)!}$ has mass one,
mean $n$, and variance $n$, as is checked by integrating by parts
the moments of orders zero, one, and two.
Fix $T>0$ and take $h=\displaystyle\frac{t}{n}$, with $0<t\leq T$. Then
\[
\int_0^\infty\left|\frac{tr}{n}-t\right|^2d\nu_n(r)
=\frac{t^2}{n}\leq\frac{T^2}{n}.
\]
The function $s\mapsto T_sx$ is uniformly continuous on
$[0,T+1]$. Given $0<\delta<1$, split the integral according to
$|\displaystyle\frac{tr}{n}-t|<\delta$ or $|\displaystyle\frac{tr}{n}-t|\geq\delta$. If $\omega_x(\delta)$
is its modulus of continuity on that interval, contractivity and
Chebyshev's inequality give
\[
\sup_{0<t\leq T}
\|(I-(\frac{t}{n})G)^{-n}x-T_tx\|
\leq\omega_x(\delta)+\frac{2\|x\|T^2}{n\delta^2}.
\]
First let $n\to\infty$, then $\delta\to0^{+}$.
For $t=0$, both operators are the identity, so
convergence is uniform on $[0,T]$.
\end{proof}

\begin{proposition}[Bounded perturbations of generators]
\label{prop:perturbacion-acotada-generador-semigrupo}
Let $G$ be the generator of a $C_0$-semigroup $(T_t)_{t\geq0}$ on a
Banach space $X$, and suppose that
\[
\|T_t\|_{\mathcal L(X)}\leq e^{\omega t},
\qquad t\geq0,
\]
for some $\omega\in\mathbb R$. If $B\in\mathcal L(X)$, then the
operator
\[
G+B,
\qquad
\mathcal D(G+B)=\mathcal D(G),
\]
generates a $C_0$-semigroup $(S_t)_{t\geq0}$. This semigroup is given by the
Dyson--Phillips series
\begin{equation}
\label{eq:serie-dyson-phillips-perturbacion-acotada}
S_t
=
\sum_{j=0}^{\infty}S_j(t),
\qquad
S_0(t):=T_t,
\qquad
S_{j+1}(t)x
:=
\int_0^tT_{t-\tau}B S_j(\tau)x\,d\tau,
\end{equation}
which converges in operator norm uniformly for $t$ in compact
intervals. Moreover,
\[
\|S_t\|_{\mathcal L(X)}
\leq
 e^{(\omega+\|B\|)t}.
\]
In particular, if $(T_t)$ is a contraction semigroup and $\mu\geq\|B\|$, then
$(e^{-\mu t}S_t)_{t\geq0}$ is a contraction $C_0$-semigroup with
generator $G+B-\mu I$.
\end{proposition}

\begin{proof}
The operators in \eqref{eq:serie-dyson-phillips-perturbacion-acotada}
are well defined by Bochner integrals with values in $X$.
Indeed, strong continuity of $T_t$ and local boundedness imply
joint continuity of $(t,z)\mapsto T_tz$. Inductively, the change
$\tau=t\theta$ in each integral proves that $t\mapsto S_j(t)x$ is
continuous for every $x\in X$, with $S_j(0)=0$ if $j\geq1$.
The induction also gives
\[
 \|S_j(t)\|\leq e^{\omega t}\frac{\|B\|^jt^j}{j!}.
\]
For the induction step, integrate
$e^{\omega(t-\tau)}\|B\|e^{\omega\tau}\|B\|^j\displaystyle\frac{\tau^j}{j!}$
over $0\leq\tau\leq t$. Therefore the series converges in operator
norm uniformly on compact intervals, defines a strongly continuous
family with $S_0=I$, and satisfies
\[
 \|S_t\|\leq e^{(\omega+\|B\|)t}.
\]
Uniform convergence permits summation of the integrals for each $x$ and gives
\[
 S_tx=T_tx+\int_0^tT_{t-\tau}BS_\tau x\,d\tau.
\]

Fix $s\geq0$. Splitting the integral for $S_{t+s}x$ into $[0,s]$
and $[s,t+s]$, the semigroup property of $T$ gives
\[
 S_{t+s}x=T_tS_sx+\int_0^tT_{t-\tau}BS_{s+\tau}x\,d\tau.
\]
The integral formula applied to $S_sx$ gives the same identity
for $S_tS_sx$, with $S_\tau S_sx$ in the integrand. Their difference
$D_t:=S_{t+s}x-S_tS_sx$ satisfies
\[
 D_t=\int_0^tT_{t-\tau}BD_\tau\,d\tau.
\]
On $[0,R]$, set $K_R:=\sup_{0\leq t\leq R}\|D_t\|$ and
$c_R:=e^{\max\{\omega,0\}R}\|B\|$. Continuity makes the first
constant finite. Iteration of the integral inequality gives,
for every $n\in\mathbb N$,
\[
 \|D_t\|\leq K_R\frac{(c_Rt)^n}{n!},\qquad 0\leq t\leq R.
\]
Letting $n\to\infty$ gives $D_t=0$. Since $R$ is arbitrary,
$S_{t+s}=S_tS_s$ and $(S_t)$ is a $C_0$-semigroup.

Denote its generator by $H$. For every $x\in X$,
\[
 \frac{1}{t}S_1(t)x
 =\frac{1}{t}\int_0^tT_{t-\tau}BT_\tau x\,d\tau\longrightarrow Bx.
\]
This convergence follows from
\[
 \sup_{0\leq\tau\leq t}\|T_{t-\tau}BT_\tau x-Bx\|
 \leq e^{\max\{\omega,0\}}\|B\|
       \sup_{0\leq\tau\leq t}\|T_\tau x-x\|
       +\sup_{0\leq r\leq t}\|(T_r-I)Bx\|\longrightarrow0
\]
for $0<t\leq1$. Moreover,
\[
 \left\|\sum_{j=2}^{\infty}S_j(t)\right\|
 \leq\frac{1}{2} e^{\max\{\omega,0\}+\|B\|}\|B\|^2t^2,
 \qquad 0<t\leq1.
\]
For $x\in\mathcal D(G)$, dividing the series by $t$ proves that
$x\in\mathcal D(H)$ and $Hx=Gx+Bx$.

Conversely, for $x\in\mathcal D(H)$ the integral formula gives
\[
 \frac{T_tx-x}{t}
 =\frac{S_tx-x}{t}-\frac{1}{t}\int_0^tT_{t-\tau}BS_\tau x\,d\tau.
\]
The first term converges to $Hx$. For the second,
\[
 \sup_{0\leq\tau\leq t}\|T_{t-\tau}BS_\tau x-Bx\|
 \leq e^{\max\{\omega,0\}}\|B\|
       \sup_{0\leq\tau\leq t}\|S_\tau x-x\|
       +\sup_{0\leq r\leq t}\|(T_r-I)Bx\|\longrightarrow0.
\]
Its average therefore converges to $Bx$. By the definition of $G$,
$x\in\mathcal D(G)$ and $Gx=Hx-Bx$. Consequently,
$\mathcal D(H)=\mathcal D(G)$ and $H=G+B$.
Finally, the difference quotient of $e^{-\mu t}S_t$ gives generator
$H-\mu I$ on the same domain; if $\omega=0$ and $\mu\geq\|B\|$,
the bound on $S_t$ shows that this semigroup is contractive.
\end{proof}

\subsection{Bessel potentials associated with a semigroup}

\begin{lemma}[Beta--gamma identity]\label{lem:identidad-beta-gamma}
If $a,b>0$, then
\[
\int_0^1\theta^{a-1}(1-\theta)^{b-1}\,d\theta
=
\frac{\Gamma(a)\Gamma(b)}{\Gamma(a+b)}.
\]
\end{lemma}

\begin{proof}
By Tonelli's theorem,
\[
\Gamma(a)\Gamma(b)
=
\int_0^\infty\int_0^\infty x^{a-1}y^{b-1}e^{-(x+y)}\,dy\,dx.
\]
With the change of variables $r=x+y$, $\theta=\displaystyle\frac{x}{x+y}$, whose Jacobian is $r$, the right-hand side becomes
\[
\left(\int_0^\infty r^{a+b-1}e^{-r}\,dr\right)
\left(\int_0^1\theta^{a-1}(1-\theta)^{b-1}\,d\theta\right).
\]
The first integral is $\Gamma(a+b)$.
\end{proof}

\begin{theorem}[Bessel potentials of a contraction semigroup]\label{teo:escala-bessel-semigrupo-contraccion}\index{Bessel potential!of a contraction semigroup}
Let $(T_t)_{t\geq0}$ be a contraction $C_0$-semigroup on a Banach space $X$, and let $-A$ be its generator. For $s>0$, define
\[
J_sx
:=
\frac{1}{\Gamma\left(\frac{s}{2}\right)}
\int_0^\infty t^{\frac{s}{2}-1}e^{-t}T_tx\,dt.
\]
Then:
\begin{enumerate}[label=(\alph*)]
\item $J_s\in\mathcal L(X)$ and $\|J_s\|\leq1$;
\item $J_sJ_r=J_{s+r}$ for $r,s>0$;
\item $J_2=(I+A)^{-1}$ and $J_{2m}=(I+A)^{-m}$ for $m\in\mathbb N$;
\item each $J_s$ is injective and has dense range in $X$.
\end{enumerate}

For $s>0$, let $X_s:=J_s(X)$ with the norm $\|J_sx\|_{X_s}:=\|x\|_X$; let $X_0:=X$ and, for $s<0$, let $X_s$ be the completion of $X$ with respect to $\|x\|_{X_s}:=\|J_{-s}x\|_X$. Then each $X_s$ is a Banach space, $X_{s_1}\hookrightarrow X_{s_0}$ continuously if $s_1>s_0$, and $J_r$ extends uniquely to an isometric isomorphism
\[
J_r\colon X_s\longrightarrow X_{s+r}
\]
for any $s\in\mathbb R$ and $r>0$.
\end{theorem}

\begin{proof}
The estimate
\[
 \|J_sx\|\leq
 \frac{\|x\|}{\Gamma(\frac{s}{2})}\int_0^\infty t^{\frac{s}{2}-1}e^{-t}\,dt
 =\|x\|
\]
proves (a). The double integral representing $J_sJ_rx$ is absolutely
integrable in norm, since $\|T_{t+\tau}x\|\leq\|x\|$.
Fubini, the change $a=t+\tau$, $\theta=\displaystyle\frac{t}{a}$, and
Lemma~\ref{lem:identidad-beta-gamma} give
\[
 J_sJ_rx
 =\frac{1}{\Gamma(\frac{s}{2})\Gamma(\frac{r}{2})}
 \int_0^\infty\int_0^\infty
 t^{\frac{s}{2}-1}\tau^{\frac{r}{2}-1}e^{-(t+\tau)}T_{t+\tau}x\,d\tau\,dt
 =J_{s+r}x.
\]
This proves (b). Proposition~\ref{prop:resolvente-transformada-laplace-semigrupo}
with $G=-A$ and $\lambda=1$ gives $J_2=(I+A)^{-1}$, and (b) implies (c).

If $J_sx=0$, choose an integer $m$ with $2m>s$. Then
$J_{2m}x=J_{2m-s}J_sx=0$. Since $(I+A)^{-m}$ is injective,
$x=0$. Moreover,
\[
 \operatorname{Ran}(J_{2m})
 =\mathcal D((I+A)^m)=\mathcal D(A^m)
 \subseteq\operatorname{Ran}(J_s).
\]
Equality of the domains follows by induction. On $\mathcal D(A^m)$,
the binomial expansion shows that $(I+A)^mx-A^mx\in\mathcal D(A)$.
Thus the additional condition $(I+A)^mx\in\mathcal D(A)$ is equivalent
to $A^mx\in\mathcal D(A)$, giving the step from $m$ to $m+1$.
Density of this domain follows from
Proposition~\ref{prop:propiedades-generador-semigrupo-C0}.
This proves (d). Also set $J_0:=I$.

For $s\geq0$, define $U_s\colon X_s\to X$ by
$U_s(J_sx)=x$. It is a surjective isometry, so $X_s$ is complete.
For $s<0$, the operator $J_{-s}$ is an isometry from $X$, equipped
with the norm defining $X_s$, into $X$. It extends to an isometry
$U_s\colon X_s\to X$. Its range is closed by completeness of
$X_s$ and contains the dense range of $J_{-s}$; hence $U_s$ is
surjective as well. We have constructed a surjective isometry
$U_s$ for every $s\in\mathbb R$.

For $a>b$, define
\[
 i_{a,b}:=U_b^{-1}J_{a-b}U_a\colon X_a\longrightarrow X_b.
\]
This operator has norm at most one and is injective because $J_{a-b}$
is injective. For $a>b\geq0$, it agrees with the usual inclusion:
$J_ax=J_bJ_{a-b}x$. If $a\geq0>b$, its value on $J_ax$ agrees
with the image of that vector in the completion $X_b$, since
$J_{-b}J_ax=J_{a-b}x$. If $0>a>b$, on the canonical subspace
$X$ the identity $J_{a-b}J_{-a}=J_{-b}$ proves the same compatibility.
Moreover, for $a>b>c$,
\[
 i_{b,c}i_{a,b}=U_c^{-1}J_{b-c}J_{a-b}U_a=i_{a,c}.
\]
These are the continuous compatible inclusions of the scale,
also when both indices are negative or have opposite signs.

For $r>0$ and $s\in\mathbb R$, define the scale shift by
\[
 J_r^{(s)}:=U_{s+r}^{-1}U_s\colon X_s\longrightarrow X_{s+r}.
\]
It is a surjective isometry with inverse $U_s^{-1}U_{s+r}$.
If $s\geq0$, it agrees with $J_r$ on $X_s$, since
$J_rJ_sx=J_{s+r}x$. If $s<0$ and $s+r\leq0$, for $x\in X$ we have
$U_{s+r}J_rx=J_{-(s+r)}J_rx=J_{-s}x=U_sx$.
If $s<0<s+r$, the identity $J_rx=J_{s+r}J_{-s}x$ again gives
$U_{s+r}J_rx=U_sx$. By density of $X$ in the negative spaces,
these shifts are the unique continuous extensions of the original
operator. The formulas involving $U_s$ also show compatibility
with the inclusions $i_{a,b}$. This proves all assertions about the scale.
\end{proof}

\begin{lemma}[Integral formula for the square root of a generator]
\label{lem:formula-integral-raiz-semigrupo}
Let $-A$ be the generator of a contraction $C_0$-semigroup
$(e^{-tA})_{t\geq0}$ on a Banach space $X$, let $a>0$, and set
$B:=aI+A$. Define the bounded operator
\[
C:=\frac{1}{\sqrt\pi}\int_0^\infty
t^{-\frac{1}{2}}e^{-tB}\,dt.
\]
Then $C^2=B^{-1}$; in particular, $C$ is injective and has dense
range. Set
\[
B^{-\frac{1}{2}}:=C,\qquad
\mathcal D(B^{\frac{1}{2}}):=\operatorname{Ran}(C),\qquad
B^{\frac{1}{2}}:=C^{-1}\colon\operatorname{Ran}(C)\longrightarrow X.
\]
For every $x\in\mathcal D(B)$, the Bochner integral
\begin{equation}
\label{eq:formula-integral-raiz-semigrupo}
\frac{1}{2\sqrt\pi}
\int_0^\infty
t^{-\frac{3}{2}}\bigl(x-e^{-tB}x\bigr)\,dt
\end{equation}
converges in $X$ and equals $B^{\displaystyle\frac{1}{2}}x$. Since
$\mathcal D(B)$ is dense in $\mathcal D(B^{\displaystyle\frac{1}{2}})$ for the norm
$\|x\|_{\displaystyle\frac{1}{2}}:=\|B^{\displaystyle\frac{1}{2}}x\|_X$, which is equivalent to the graph norm, this
identity determines $B^{\displaystyle\frac{1}{2}}$ by closure.
\end{lemma}

\begin{proof}
If $x\in\mathcal D(B)$, then
\[
x-e^{-tB}x
=
\int_0^t e^{-sB}Bx\,ds.
\]
Thus
\[
\|x-e^{-tB}x\|_X
\leq
\min\{t\|Bx\|_X,2\|x\|_X\}.
\]
This bound proves convergence of
\eqref{eq:formula-integral-raiz-semigrupo} both near zero and at
infinity.

The bound $\|e^{-tB}\|\leq e^{-at}$ allows Fubini's theorem to be applied to the integral
defining $C^2$. With $r=t+s$ and $\theta=\displaystyle\frac{t}{r}$, the beta identity gives
\begin{align*}
C^2x
&=\frac{1}{\pi}\int_0^\infty\int_0^\infty
t^{-\frac{1}{2}}s^{-\frac{1}{2}}e^{-(t+s)B}x\,ds\,dt\\
&=\frac{1}{\pi}\int_0^\infty e^{-rB}x
\left(\int_0^1\theta^{-\frac{1}{2}}(1-\theta)^{-\frac{1}{2}}\,d\theta\right)dr\\
&=\int_0^\infty e^{-rB}x\,dr=B^{-1}x.
\end{align*}
Hence $C$ is injective and
$\operatorname{Ran}(C)\supset\operatorname{Ran}(C^2)
=\operatorname{Ran}(B^{-1})=\mathcal D(B)$ is dense.

Let $f\in X$ and take $x=B^{-1}f$. Fubini's theorem for Bochner
integrals and the preceding identity give
\begin{align*}
&\frac{1}{2\sqrt\pi}
\int_0^\infty
t^{-\frac{3}{2}}\bigl(B^{-1}f-e^{-tB}B^{-1}f\bigr)\,dt\\
&\qquad=
\frac{1}{2\sqrt\pi}
\int_0^\infty t^{-\frac{3}{2}}
\int_0^t e^{-sB}f\,ds\,dt\\
&\qquad=
\frac{1}{\sqrt\pi}
\int_0^\infty s^{-\frac{1}{2}}e^{-sB}f\,ds
=
B^{-\frac{1}{2}}f.
\end{align*}
In the last equality we used the definition of $C$. Since
$B^{-1}f=C^2f=C(Cf)$, the computed value is
$Cf=C^{-1}x=B^{\displaystyle\frac{1}{2}}x$. Finally, under the isometry
\[
C^{-1}\colon
(\mathcal D(B^{\frac{1}{2}}),\|\cdot\|_{\frac{1}{2}})\longrightarrow X,
\]
the subspace $\mathcal D(B)=\operatorname{Ran}(C^2)$ corresponds to
$\operatorname{Ran}(C)$, which is dense in $X$. This proves the asserted
density. Moreover,
$\|x\|_X\leq\|C\|_{\mathcal L(X)}\|x\|_{\displaystyle\frac{1}{2}}$, so
$\|\cdot\|_{\displaystyle\frac{1}{2}}$ is equivalent to the graph norm.
\end{proof}

\paragraph{Uniqueness of the Laplace transform for finite measures.}
If a finite signed or complex Borel measure $\nu$ on $[0,\infty)$ satisfies
$\displaystyle\int_{[0,\infty)} e^{-\lambda t}\,d\nu(t)=0$ for every $\lambda>0$, then $\nu=0$.
Indeed, dominated convergence as $\lambda\to0^{+}$ also gives
$\nu([0,\infty))=0$. Push $\nu$ forward by $t\mapsto e^{-t}$ to $(0,1]$
and extend it to a measure $\widetilde\nu$ on $[0,1]$ with zero mass at $0$.
For every integer $k\geq0$, its moment $\displaystyle\int_{[0,1]}x^k\,d\widetilde\nu(x)$ is
zero. Polynomials are uniformly dense in $C([0,1])$ by the Weierstrass
approximation theorem, so $\widetilde\nu$ annihilates all continuous
functions. Finite Borel measures on the compact metric interval are
regular; uniqueness in the Riesz representation theorem
(Theorem~\ref{teo: representacion de Riesz en Cc}), applied to real and
imaginary parts, gives $\widetilde\nu=0$. The map to $(0,1]$ is a Borel
isomorphism, hence $\nu=0$. Applying this argument to a difference proves
that equal Laplace transforms have equal representing finite measures.

\begin{proposition}[Laplace calculus for a contraction semigroup]
\label{prop:calculo-laplace-medidas-semigrupo}
Let $(T_t)_{t\geq0}$ be a contraction $C_0$-semigroup on a
real or complex Banach space $X$. Let $\nu$ be a finite Borel measure
on $[0,\infty)$, signed if $X$ is real, and signed or complex
if $X$ is complex. Then the Bochner integral
\[
\mathcal L_\nu(T)x
:=
\int_{[0,\infty)}T_tx\,d\nu(t)
\]
defines a bounded operator on $X$, and
\[
\|\mathcal L_\nu(T)\|_{\mathcal L(X)}
\leq
|\nu|([0,\infty)).
\]
If $\mu$ is another finite measure with values in the scalar field of $X$ and $\nu*\mu$ denotes their convolution, then
\[
\mathcal L_\nu(T)\mathcal L_\mu(T)
=
\mathcal L_{\nu*\mu}(T).
\]
In particular, the product of two Laplace transforms of finite
measures is represented by the product of the corresponding
operators.
\end{proposition}

\begin{proof}
For a simple function $g=\sum_{j=1}^Nx_j\mathbf1_{E_j}$,
with disjoint $E_j$, integration with respect to $\nu$ is
$\sum_{j=1}^Nx_j\nu(E_j)$ and satisfies
\[
 \left\|\int_{[0,\infty)}g\,d\nu\right\|_X
 \leq\sum_{j=1}^N\|x_j\|_X|\nu(E_j)|
 \leq\int_{[0,\infty)}\|g\|_X\,d|\nu|.
\]
The scalars $\nu(E_j)$ belong to the field of $X$ by the
hypothesis in the statement. Approximation by simple functions
in $L^1(|\nu|;X)$ extends this integral uniquely to
Bochner integrable functions and preserves the bound.
For a signed measure it agrees with the difference of the
integrals with respect to $\nu^+$ and $\nu^-$; for a complex
measure it applies to its real and imaginary parts.

Strong continuity of $t\mapsto T_tx$ ensures Bochner
measurability, and contractivity gives
\[
\int_{[0,\infty)}\|T_tx\|_X\,d|\nu|(t)
\leq
|\nu|([0,\infty))\|x\|_X.
\]
This proves existence of the integral and the bound. For the product
rule, Fubini's theorem for Bochner integrals gives
\[
\begin{split}
\mathcal L_\nu(T)\mathcal L_\mu(T)x
&=
\int_{[0,\infty)}\int_{[0,\infty)}T_{t+s}x\,d\mu(s)\,d\nu(t)\\
&=
\int_{[0,\infty)}T_rx\,d(\nu*\mu)(r),
\end{split}
\]
where we used $T_tT_s=T_{t+s}$. Scalar Fubini shows that the Laplace
transform of $\nu*\mu$ is the product of the two transforms. The uniqueness
just proved identifies any other finite measure representing that same
product with $\nu*\mu$, which proves the last assertion also at the level
of scalar functions.
\end{proof}

\begin{lemma}[Independence of the positive shift]
\label{lem:independencia-desplazamiento-bessel-Lp}
Let $-A$ be the generator of a contraction $C_0$-semigroup
$(e^{-tA})_{t\geq0}$ on a Banach space $X$. For $a>0$ and $s>0$,
define
\[
J_s^{(a)}x
:=
\frac{1}{\Gamma\left(\frac{s}{2}\right)}
\int_0^\infty
t^{\frac{s}{2}-1}e^{-at}e^{-tA}x\,dt.
\]
Each $J_s^{(a)}$ is injective. Equip its range with the norm
\[
\bigl\|J_s^{(a)}x\bigr\|_{\operatorname{Ran}(J_s^{(a)})}
:=
\|x\|_X.
\]
If $a,b>0$, then
\[
\operatorname{Ran}(J_s^{(a)})
=
\operatorname{Ran}(J_s^{(b)}),
\]
and the two range norms are equivalent.
\end{lemma}

\begin{proof}
Set $\alpha=\displaystyle\frac{s}{2}$. The Laplace formula for negative powers
gives
\[
J_s^{(a)}=(aI+A)^{-\alpha}.
\]
The semigroup property of these powers and the identity
$J_{2m}^{(a)}=(aI+A)^{-m}$ show that $J_s^{(a)}$ is injective. Indeed,
if $J_s^{(a)}x=0$, choose $m\in\mathbb N$ with $2m>s$ and obtain
\[
(aI+A)^{-m}x
=
J_{2m-s}^{(a)}J_s^{(a)}x
=0.
\]
The resolvent of $-A$ is injective, so $x=0$. Thus the norm
in the statement is well defined.
Begin with two shifts $c,d>0$ satisfying
$|d-c|<c$. For $\lambda\geq0$, the binomial series converges absolutely, and
\[
\left(\frac{d+\lambda}{c+\lambda}\right)^\alpha
=
\left(1+\frac{d-c}{c+\lambda}\right)^\alpha
=
1+
\sum_{k=1}^\infty
\binom{\alpha}{k}(d-c)^k(c+\lambda)^{-k}.
\]
For $k\geq1$, we have
\[
(c+\lambda)^{-k}
=
\int_0^\infty e^{-\lambda t}
\frac{t^{k-1}e^{-ct}}{(k-1)!}\,dt.
\]
Thus the preceding function is the Laplace transform of
\[
\nu_{c,d}
:=
\delta_0+
\sum_{k=1}^\infty
\binom{\alpha}{k}(d-c)^k
\frac{t^{k-1}e^{-ct}}{(k-1)!}\,dt.
\]
This series converges in total variation because
\[
\sum_{k=1}^\infty
\left|\binom{\alpha}{k}\right|
\left(\frac{|d-c|}{c}\right)^k
<\infty.
\]
Proposition~\ref{prop:calculo-laplace-medidas-semigrupo} then defines
a bounded operator $Q_{c,d}$ whose scalar multiplier is
$\left(\displaystyle\frac{d+\lambda}{c+\lambda}\right)^\alpha$. The product rule gives
\begin{equation}
\label{eq:factorizacion-desplazamientos-apendice}
(cI+A)^{-\alpha}
=
(dI+A)^{-\alpha}Q_{c,d}.
\end{equation}

For arbitrary shifts $a,b>0$, choose $N\in\mathbb N$ so
that $\left(\displaystyle\frac{b}{a}\right)^{\displaystyle\frac{1}{N}}<2$, and set
\[
a_j:=a\left(\frac{b}{a}\right)^{\frac{j}{N}},
\qquad j\in\{0,\dots,N\}.
\]
Then $|a_{j+1}-a_j|<a_j$ for every $j$. Iterating
\eqref{eq:factorizacion-desplazamientos-apendice} gives a bounded
operator $Q_{a,b}$ such that
\[
(aI+A)^{-\alpha}
=
(bI+A)^{-\alpha}Q_{a,b}.
\]
If $u=J_s^{(a)}x$, this identity writes
$u=J_s^{(b)}Q_{a,b}x$ and, consequently,
\[
\|u\|_{\operatorname{Ran}(J_s^{(b)})}
\leq
\|Q_{a,b}\|\,
\|u\|_{\operatorname{Ran}(J_s^{(a)})}.
\]
Interchanging $a$ and $b$ gives the reverse inclusion and estimate.
\end{proof}

\begin{proposition}[Complex-linear duality of consistent scales]
\label{prop:dualidad-escala-bessel-consistente}
Let $X$ and $Y$ be complex Banach spaces, and let
$\mathfrak b\colon X\times Y\to\mathbb C$ be a continuous bilinear form such
that
\[
 Y\longrightarrow X',\qquad w\longmapsto\mathfrak b(\,\cdot\,,w),
\]
is a complex-linear isometric isomorphism. Let
$(T_t)_{t\geq0}$ and $(S_t)_{t\geq0}$ be contraction semigroups on $X$ and
$Y$, respectively, and suppose that
\begin{equation}
\label{eq:semigrupos-duales-conjugados-apendice}
 \mathfrak b(T_tu,w)=\mathfrak b(u,S_tw)
\end{equation}
for $t\geq0$ and $u,w$ in respective dense invariant subspaces. If
$X_s$ and $Y_s$ are the scales from
Theorem~\ref{teo:escala-bessel-semigrupo-contraccion}, then the duality
$\mathfrak b$ extends uniquely to a complex-linear isometric
isomorphism
\[
 Y_{-s}\xrightarrow{\ \cong\ }(X_s)',
 \qquad s\in\mathbb R.
\]
In the application to a Hermitian bundle $\mathbf{E}$, take
$X=L^p(M,\mathbf{E})$, $Y=L^{p'}(M,\overline{\mathbf{E}})$, and
$\displaystyle \mathfrak b(u,\overline v)=\int_M\mathbf{h}_{\mathbf{E}}(u,v)\,d\lambda_{\mathbf{g}}$; thus the
space representing the complex dual of a scale on $\mathbf{E}$ is a
scale on $\overline{\mathbf{E}}$, not another complex-linear copy of the scale
on $\mathbf{E}$.
\end{proposition}

\begin{proof}
The identity
\eqref{eq:semigrupos-duales-conjugados-apendice}, first on the dense
subspaces and then by continuity, and the Bochner integral defining the
Bessel potentials imply, for every $r>0$,
\begin{equation}
\label{eq:potenciales-bessel-duales-conjugados}
 \mathfrak b(J_r^Xu,w)=\mathfrak b(u,J_r^Yw).
\end{equation}

First suppose that $s>0$. If $w\in Y$ and $v=J_s^Xu\in X_s$, then
\[
 |\mathfrak b(v,w)|
 =|\mathfrak b(u,J_s^Yw)|
 \leq\|v\|_{X_s}\|w\|_{Y_{-s}}.
\]
By density, the map $w\mapsto\mathfrak b(\,\cdot\,,w)$ extends
isometrically from $Y_{-s}$ to $(X_s)'$. To see that it is surjective, let
$F\in(X_s)'$. The functional $u\mapsto F(J_s^Xu)$ on $X$ is represented
by a unique $z\in Y$. Since
$J_s^Y\colon Y_{-s}\to Y$ is an isometric isomorphism, there exists a unique
$w\in Y_{-s}$ with $J_s^Yw=z$. Identity
\eqref{eq:potenciales-bessel-duales-conjugados}, extended by density,
gives $F(v)=\mathfrak b(v,w)$ and preserves the norm.

Now let $s=-r<0$. If $w=J_r^Yz\in Y_r$, then, for $u\in X$,
\[
 |\mathfrak b(u,w)|
 =|\mathfrak b(J_r^Xu,z)|
 \leq\|u\|_{X_{-r}}\|z\|_Y
 =\|u\|_{X_s}\|w\|_{Y_r}.
\]
This extends $w$ to an element of $(X_s)'$. Conversely, if
$F\in(X_{-r})'$, its restriction to $X$ is represented by some
$w\in Y$. The bound on $F$ makes the functional $J_r^Xu\mapsto\mathfrak b(u,w)$ continuous on $\operatorname{Ran}(J_r^X)$;
since that range is dense in
$X$, it extends to a functional represented by some $z\in Y$. Then
\eqref{eq:potenciales-bessel-duales-conjugados} gives
$\mathfrak b(u,w)=\mathfrak b(u,J_r^Yz)$ for every $u\in X$, hence
$w=J_r^Yz\in Y_r$. The argument preserves norms. The case $s=0$ is the
duality $Y\cong X'$ already mentioned.
\end{proof}

\section{Measure spaces}\label{ap:espacios-de-medida}

We now turn to the basic concepts of measure theory. The aim of this section is to establish a unified language for defining Lebesgue and Bochner integrals, $L^{p}(X)$ spaces, and measures induced by geometric structures.

\begin{definition}\label{def:b4-espacios-medida-algebra-sobre}\index{sigma algebra@sigma-algebra}
Let $X$ be a set. A family $\mathcal{A}\subseteq\mathcal{P}(X)$ is a \textbf{$\sigma$-algebra on $X$} if $X\in\mathcal{A}$, if $X\setminus A\in\mathcal{A}$ for each $A\in\mathcal{A}$, and if $\displaystyle\bigcup_{n=1}^{\infty}A_{n}\in\mathcal{A}$ for every sequence $(A_{n})_{n\in\mathbb{N}}\subseteq\mathcal{A}$. In this case, $(X,\mathcal{A})$ is called a \textbf{measurable space}, and the elements of $\mathcal{A}$ are called \textbf{$\mathcal{A}$-measurable} sets, or simply \textbf{measurable} sets when the $\sigma$-algebra is understood.
\end{definition}

\begin{definition}\label{def:b4-espacios-medida-medida-en}\index{measure}
Let $(X,\mathcal{A})$ be a measurable space. A \textbf{measure} on $(X,\mathcal{A})$ is a function $\mu\colon \mathcal{A}\longrightarrow[0,\infty]$ such that $\mu(\varnothing)=0$ and $\displaystyle\mu\left(\bigcup_{n=1}^{\infty}A_{n}\right)=\displaystyle\sum_{n=1}^{\infty}\mu(A_{n})$ for every sequence $(A_{n})_{n\in\mathbb{N}}\subseteq\mathcal{A}$ of pairwise disjoint sets. The triple $(X,\mathcal{A},\mu)$ is called a \textbf{measure space}.
\end{definition}

Before constructing specific measures, we recall how the $\sigma$-algebras determining measurable sets are generated.

\begin{definition}\label{def:b4-espacios-medida-algebra-generada-por}\index{generated sigma algebra@generated sigma-algebra}\index{Borel sigma algebra@Borel $\sigma$-algebra}
Let $X$ be a set, and let $\mathcal{S}\subseteq\mathcal{P}(X)$. The \textbf{$\sigma$-algebra generated by $\mathcal{S}$} is $\displaystyle\sigma(\mathcal{S}):=\bigcap\{\mathcal{A}\subseteq\mathcal{P}(X)\mid \mathcal{S}\subseteq\mathcal{A}\text{ and $\mathcal{A}$ is a $\sigma$-algebra on $X$}\}$. If $(X,\tau)$ is a topological space, the \textbf{Borel $\sigma$-algebra} of $X$ is $\mathcal{B}(X):=\sigma(\tau)$, and its elements are called Borel sets.
\end{definition}

\begin{definition}\label{def:b4-espacios-medida-nulo}\index{set of measure zero}\index{complete measure}
Let $(X,\mathcal{A},\mu)$ be a measure space. A set $N\in\mathcal{A}$ is \textbf{$\mu$-null} if $\mu(N)=0$. The measure $\mu$, or equivalently the space $(X,\mathcal{A},\mu)$, is \textbf{complete} if every subset of a $\mu$-null set belongs to $\mathcal{A}$; by monotonicity, that subset is then also $\mu$-null.
\end{definition}

\begin{definition}\label{def:b4-espacios-medida-medida}\index{measure!finite}\index{measure!sigma-finite@$\sigma$-finite}
Let $(X,\mathcal{A},\mu)$ be a measure space. We say that $\mu$ is \textbf{finite} if $\mu(X)<\infty$, and \textbf{$\sigma$-finite} if there exists a sequence $(A_{n})_{n\in\mathbb{N}}\subseteq\mathcal{A}$ such that $X=\displaystyle\bigcup_{n=1}^{\infty}A_{n}$ and $\mu(A_{n})<\infty$ for every $n\in\mathbb{N}$.
\end{definition}

Carathéodory's construction produces a complete measure from a function initially defined on all subsets of $X$.
Complete proofs of this construction, regularity of Lebesgue measure, and the signed measure decomposition theorems stated below can be found in \cite{Grabinsky,Bogachev2007,Cohn2013}.

\begin{definition}\label{def:b4-espacios-medida-medida-exterior-en}\index{outer measure}
An \textbf{outer measure} on a set $X$ is a function $\mu^{*}\colon \mathcal{P}(X)\longrightarrow[0,\infty]$ such that $\mu^{*}(\varnothing)=0$, $\mu^{*}(A)\leq\mu^{*}(B)$ whenever $A\subseteq B$, and $\displaystyle\mu^{*}\left(\bigcup_{k=1}^{\infty}A_{k}\right)\leq\displaystyle\sum_{k=1}^{\infty}\mu^{*}(A_{k})$ for every sequence $(A_{k})_{k\in\mathbb{N}}\subseteq\mathcal{P}(X)$.
\end{definition}

\begin{theorem}\label{teo:b4-espacios-medida-medida-exterior-de-lebesgue-dimensional}\index{Lebesgue outer measure dimensional@$n$-dimensional Lebesgue outer measure}
For a compact rectangle $R=\displaystyle\prod_{i=1}^{n}[a_{i},b_{i}]$, define $\displaystyle\operatorname{vol}_{n}(R):=\prod_{i=1}^{n}(b_{i}-a_{i})$. Given $A\subseteq\mathbb{R}^{n}$, let $\mathcal{C}(A)$ be the family of sequences $(R_{j})_{j\in\mathbb{N}}$ of compact rectangles such that $A\subseteq\displaystyle\bigcup_{j=1}^{\infty}R_{j}$, and define $\displaystyle\lambda_{n}^{*}(A):=\inf_{(R_{j})\in\mathcal{C}(A)}\displaystyle\sum_{j=1}^{\infty}\operatorname{vol}_{n}(R_{j})$. Then $\lambda_{n}^{*}$ is an outer measure on $\mathbb{R}^{n}$, called \textbf{$n$-dimensional Lebesgue outer measure}.
\end{theorem}

Not every subset of a space equipped with an outer measure is necessarily measurable. The following condition selects precisely those on which the outer measure becomes countably additive.

\begin{definition}\label{def:b4-espacios-medida-medible}\index{measurable set!in the sense of Carathéodory}
Let $\mu^{*}$ be an outer measure on $X$. A subset $A\subseteq X$ is \textbf{$\mu^{*}$-measurable} if $\mu^{*}(D)=\mu^{*}(D\cap A)+\mu^{*}(D\setminus A)$ for each $D\subseteq X$. Denote by $\mathcal{A}(\mu^{*})$ the family of all $\mu^{*}$-measurable subsets of $X$.
\end{definition}

\begin{theorem}[Carathéodory]\label{teo:b4-espacios-medida-medida-exterior-algebra-medida-completa-dicha}\index{outer measure!induced complete measure}
Let $\mu^{*}$ be an outer measure on $X$. Then $\mathcal{A}(\mu^{*})$ is a $\sigma$-algebra on $X$, and the restriction $\mu^{*}\restriction_{\mathcal{A}(\mu^{*})}$ is a complete measure. This restriction is called the \textbf{measure induced by $\mu^{*}$}.
\end{theorem}

\begin{definition}\label{def:b4-espacios-medida-medida-de-lebesgue}\index{Lebesgue measure}\index{Lebesgue measurable set}
The $\sigma$-algebra of \textbf{Lebesgue measurable} subsets of $\mathbb{R}^{n}$ is $\mathcal{L}(\mathbb{R}^{n}):=\mathcal{A}(\lambda_{n}^{*})$, and \textbf{$n$-dimensional Lebesgue measure} is the complete measure $\lambda_{n}:=\lambda_{n}^{*}\restriction_{\mathcal{L}(\mathbb{R}^{n})}$. When there is no risk of confusion, we simply write $dx$ in place of $d\lambda_{n}(x)$.
\end{definition}

To relate measure theory to topology, we introduce regularity and local finiteness.

\begin{definition}\label{def:b4-espacios-medida-regulares}\index{regular measure}\index{locally finite measure}
Let $X$ be a topological space, and let $(X,\mathcal{A},\mu)$ be a measure space such that $\mathcal{B}(X)\subseteq\mathcal{A}$. The measure $\mu$ is \textbf{regular} if, for each $A\in\mathcal{A}$, we have $\displaystyle\mu(A)=\inf\{\mu(U)\mid U\subseteq X\text{ is open and }A\subseteq U\}=\displaystyle\sup\{\mu(K)\mid K\subseteq A\text{ is compact}\}$. The measure $\mu$ is \textbf{locally finite} if each $x\in X$ has an open neighborhood $U$ such that $\mu(U)<\infty$.
\end{definition}

\begin{proposition}\label{medida localmente finita}
Let $X$ be a topological space, let $(X,\mathcal{A},\mu)$ be a measure space, and suppose that $\mathcal{B}(X)\subseteq\mathcal{A}$. If $\mu$ is locally finite, then $\mu(K)<\infty$ for every compact set $K\subseteq X$. Conversely, if each point of $X$ has a compact neighborhood and $\mu(K)<\infty$ for every compact set $K\subseteq X$, then $\mu$ is locally finite.
\end{proposition}
\begin{proof}
If $K\subseteq X$ is compact, for each $x\in K$ there exists an open neighborhood $U_{x}$ with $\mu(U_{x})<\infty$. By compactness, there exist $x_{1},\dots,x_{r}\in K$ such that $K\subseteq\displaystyle\bigcup_{i=1}^{r}U_{x_{i}}$, so $\displaystyle\mu(K)\leq\displaystyle\sum_{i=1}^{r}\mu(U_{x_{i}})<\infty$. Conversely, let $x\in X$, and let $K$ be a compact neighborhood of $x$. There exists an open set $U$ such that $x\in U\subseteq K$, and then $\mu(U)\leq\mu(K)<\infty$.
\end{proof}

\begin{definition}\label{def:b4-espacios-medida-casi-en-todas-partes-en}\index{almost everywhere}
Let $(X,\mathcal{A},\mu)$ be a measure space, and let $P(x)$ be a property defined for $x\in X$. We say that $P$ holds \textbf{$\mu$-almost everywhere}, or for \textbf{$\mu$-almost every} $x\in X$, if there exists a set $N\in\mathcal{A}$ with $\mu(N)=0$ such that $P(x)$ holds for every $x\in X\setminus N$. When the measure is understood, we simply write \textbf{almost everywhere}.
\end{definition}

\begin{theorem}\label{teo:b4-espacios-medida-medida-lebesgue-regular}\index{Lebesgue measure!regularity}
Lebesgue measure $\lambda_{n}$ is regular and locally finite.
\end{theorem}

Regularity of $\lambda_{n}$ allows any Lebesgue measurable set to be approximated by Borel sets, changing it only on a null set. Recall that an $F_{\sigma}$ set is a countable union of closed sets, and a $G_{\delta}$ set is a countable intersection of open sets.

\begin{theorem}\label{teorema 5.5}\index{Lebesgue measure!approximation by $F_\sigma$ and $G_\delta$ sets}
Let $A\in\mathcal{L}(\mathbb{R}^{n})$. There exist an $F_{\sigma}$ set, denoted by $F$, and a $G_{\delta}$ set, denoted by $G$, such that $F\subseteq A\subseteq G$ and $\lambda_{n}(A\setminus F)=\lambda_{n}(G\setminus A)=0$. In particular, $\lambda_{n}(F)=\lambda_{n}(A)=\lambda_{n}(G)$. If $A$ is bounded, $G$ may be chosen bounded.
\end{theorem}

\begin{definition}\label{def:b4-espacios-medida-localmente-compacto}\index{locally compact}\index{sigma compact@$\sigma$-compact}\index{Radon measure}
Let $X$ be a topological space.
\begin{enumerate}[label=(\alph*)]
\item We say that $X$ is \textbf{locally compact} if each point of $X$ has a compact neighborhood.
\item We say that $X$ is \textbf{$\sigma$-compact} if there exists a sequence $(K_{j})_{j\in\mathbb{N}}$ of compact subsets of $X$ such that $X=\displaystyle\bigcup_{j=1}^{\infty}K_{j}$.
\item If $X$ is locally compact and Hausdorff, a measure $\mu$ defined on a $\sigma$-algebra $\mathcal{A}$ with $\mathcal{B}(X)\subseteq\mathcal{A}$ is called a \textbf{Radon measure} if it is regular and locally finite. We say that $\mu$ is \textbf{massive} if, in addition, it is complete and $\mu(U)>0$ for every nonempty open set $U\subseteq X$.
\end{enumerate}
\end{definition}

The following two properties of Lebesgue measurability will be used when formulating the change-of-variables theorem on subsets that need not be open.

\begin{proposition}\label{imagen de lebesgue medible bajo c1}
Let $\Omega\subseteq\mathbb{R}^{n}$ be open, $A\in\mathcal{L}(\mathbb{R}^{n})$ with $A\subseteq\Omega$, and $f\in C^{1}(\Omega,\mathbb{R}^{m})$, where $m\geq n$. Then $f(A)\in\mathcal{L}(\mathbb{R}^{m})$.
\end{proposition}
\begin{proof}
By Theorem~\ref{teorema 5.5}, there exists an $F_{\sigma}$ set such that $F\subseteq A$ and $N:=A\setminus F$ is $\lambda_{n}$-null. Writing $F=\displaystyle\bigcup_{j=1}^{\infty}C_{j}$ with each $C_{j}$ closed and observing that $C_{j}\subseteq\Omega$, we obtain $\displaystyle F=\bigcup_{j,k=1}^{\infty}(C_{j}\cap\overline{B}_{\mathrm{euc}}(0,k))$ as a countable union of compact sets. By continuity, $f(F)$ is an $F_{\sigma}$ set and, in particular, is Borel. On the other hand, $\Omega$ admits a countable cover by convex open balls $(U_{r})_{r\in\mathbb{N}}$ with $\overline{U_{r}}\subseteq\Omega$, and the restriction of $f$ to each $U_{r}$ is Lipschitz, since $Df$ is bounded on $\overline{U_{r}}$. A Lipschitz map from a subset of $\mathbb{R}^{n}$ to $\mathbb{R}^{m}$, with $m\geq n$, sends $\lambda_{n}$-null sets to $\lambda_{m}$-null sets; to check this, cover the null set by cubes of arbitrarily small diameters and use the fact that the image of each cube is contained in a ball whose diameter is bounded by the Lipschitz constant times the diameter of the cube. Consequently, each $f(N\cap U_{r})$ is $\lambda_{m}$-null, and so is $\displaystyle f(N)=\bigcup_{r=1}^{\infty}f(N\cap U_{r})$. Finally, $f(A)=f(F)\cup f(N)$ is Lebesgue measurable.
\end{proof}

\begin{proposition}\label{lebesgue es local}
A subset $A\subseteq\mathbb{R}^{n}$ is Lebesgue measurable if and only if for each $x\in A$ there exists an open set $U_{x}\subseteq\mathbb{R}^{n}$ such that $x\in U_{x}$ and $A\cap U_{x}$ is Lebesgue measurable.
\end{proposition}
\begin{proof}
Necessity is immediate upon taking $U_{x}=\mathbb{R}^{n}$. Conversely, the family $(U_{x})_{x\in A}$ is an open cover of $A$. Since $\mathbb{R}^{n}$ is second countable and hence Lindelöf, there exists a countable subcover $(U_{x_{j}})_{j\in\mathbb{N}}$. Thus $\displaystyle A=\bigcup_{j=1}^{\infty}(A\cap U_{x_{j}})$ is Lebesgue measurable.
\end{proof}

\subsection{Signed measures}\label{ap:medidas-con-signo}

To represent linear functionals that need not be positive, we must allow measures taking values of both signs. We combine the finitely additive and countably additive versions in a single definition.

\begin{definition}\label{def:b4-espacios-medida-medida-conjunto-no-vacio-algebra}\label{def:b4-espacios-medida-medida-2}\index{finitely additive signed measure}\index{signed measure}
Let $X$ be a nonempty set, and let $\mathcal{A}$ be an algebra of subsets of $X$, that is, a nonempty family closed under complements relative to $X$ and finite unions. A function $\lambda\colon \mathcal{A}\longrightarrow[-\infty,\infty]$ is a \textbf{finitely additive signed measure} if $\lambda(\varnothing)=0$, it takes at most one of the values $\infty$ and $-\infty$, and $\lambda(E\cup F)=\lambda(E)+\lambda(F)$ whenever $E,F\in\mathcal{A}$ are disjoint. If $\mathcal{A}$ is a $\sigma$-algebra and, for every sequence $(E_{n})_{n\in\mathbb{N}}\subseteq\mathcal{A}$ of pairwise disjoint sets, $\displaystyle\lambda\left(\bigcup_{n=1}^{\infty}E_{n}\right)=\displaystyle\sum_{n=1}^{\infty}\lambda(E_{n})$ holds, then $\lambda$ is called a \textbf{signed measure}.
\end{definition}

\begin{definition}\label{def:b4-espacios-medida-medidas-mutuamente-singulares}\index{mutually singular measures}
Two positive measures $\mu$ and $\nu$ on a measurable space $(X,\mathcal{A})$ are \textbf{mutually singular}, written $\mu\perp\nu$, if there exists $P\in\mathcal{A}$ such that $\mu(X\setminus P)=0$ and $\nu(P)=0$.
\end{definition}

The Hahn--Jordan theorem shows that every signed measure decomposes canonically as the difference of two mutually singular positive measures.

\begin{theorem}[Hahn--Jordan decomposition]\label{prop:b4-espacios-medida-espacio-medible-medida-signo-medidas-positivas}\label{teo: descomposicion de Hahn--Jordan}\index{Hahn--Jordan decomposition}\index{total variation}
Let $(X,\mathcal{A})$ be a measurable space, and let $\lambda\colon \mathcal{A}\longrightarrow[-\infty,\infty]$ be a signed measure. There exist unique positive measures $\lambda^{+}$ and $\lambda^{-}$ such that $\lambda^{+}\perp\lambda^{-}$ and $\lambda=\lambda^{+}-\lambda^{-}$; moreover, one of them is finite. For each $E\in\mathcal{A}$, we have $\displaystyle\lambda^{+}(E)=\displaystyle\sup\{\lambda(F)\mid F\in\mathcal{A},\ F\subseteq E\}$ and $\displaystyle\lambda^{-}(E)=-\inf\{\lambda(F)\mid F\in\mathcal{A},\ F\subseteq E\}$. The measures $\lambda^{+}$ and $\lambda^{-}$ are called the \textbf{positive variation} and \textbf{negative variation} of $\lambda$, respectively, while $|\lambda|:=\lambda^{+}+\lambda^{-}$ is called the \textbf{total variation}. We say that $\lambda$ is finite if $|\lambda|(X)<\infty$.
\end{theorem}

\begin{definition}\label{def:b4-espacios-medida-radon-con-signo}\index{signed Radon measure}
Let $X$ be a locally compact Hausdorff space. A signed measure $\nu$ on $\mathcal{B}(X)$ is called a \textbf{signed Radon measure} if its variations $\nu^{+}$ and $\nu^{-}$ are positive Radon measures.

A \textbf{local signed Radon measure}\index{local signed Radon measure}
is a family $(\nu_U)_U$, indexed by open sets $U\subseteq X$ with
compact closure, such that each $\nu_U$ is a finite signed Radon measure
on $U$ and
\[
\nu_V\restriction_U=\nu_U
\qquad\text{if }U\subseteq V.
\]
For $u\in C_c(X)$, choose one of these open sets $U$ containing
$\supp(u)$ and define
\[
\int_X u\,d\nu:=\int_U u\,d\nu_U.
\]
The value does not depend on $U$: two choices may be compared on their union,
which also has compact closure. The same compatibility and uniqueness
in Hahn--Jordan allow the variations $\nu^+$, $\nu^-$, and $|\nu|$
to be defined on relatively compact Borel sets through the variations of
$\nu_U$. In particular, $|\nu|(K)<\infty$ for each compact set $K\subseteq X$.
A global signed Radon measure determines such a family
by restriction, but a local family need not admit values
on all Borel subsets of $X$.
\end{definition}

In what follows, $C_{c}(X)$ denotes the vector space of continuous
functions $u\colon X\longrightarrow\mathbb{R}$ with compact support.
On a general locally compact space, the regularity of representing
measures should be specified separately. When $X$ is
$\sigma$-compact, the positive measures in the following theorem are Radon
measures under the regularity convention fixed above.

\begin{theorem}[Riesz representation theorem on $C_{c}$]\label{teo: representacion de Riesz en Cc}\index{Riesz representation theorem}
Let $X$ be a locally compact Hausdorff topological space, and let $L\colon C_{c}(X)\longrightarrow\mathbb{R}$ be a linear functional.
\begin{enumerate}[label=(\alph*)]
\item If $L$ is positive, that is, $L(u)\geq0$ for every $u\in C_{c}(X)$ with $u\geq0$, then there exists a unique positive Borel measure $\mu$, finite on compact sets, outer regular on Borel sets, and inner regular on open sets and on Borel sets of finite measure, such that $\displaystyle L(u)=\int_{X}u\,d\mu$ for every $u\in C_{c}(X)$. Uniqueness is understood among measures with these regularity properties.
\item If $L$ is locally bounded, that is, if for each compact set $K\subseteq X$ there exists $c_{K}>0$ such that $|L(u)|\leq c_{K}\|u\|_{K,\infty}$ for every $u\in C_{c}(X)$ with $\operatorname{supp}(u)\subseteq K$, where $\displaystyle\|u\|_{K,\infty}:=\displaystyle\sup_{x\in K}|u(x)|$, then there exists a unique local signed Radon measure $\nu$ such that $\displaystyle L(u)=\int_{X}u\,d\nu$ for every $u\in C_{c}(X)$. Moreover, there exist positive Borel measures $\mu_{1}$ and $\mu_{2}$, with the finiteness and regularity properties in part (a), such that
\[
L(u)=\int_Xu\,d\mu_1-\int_Xu\,d\mu_2,
\qquad u\in C_c(X).
\]
On each open set $U$ with compact closure, we have
$\nu_U=\mu_1\restriction_U-\mu_2\restriction_U$.
\end{enumerate}
If $X$ is $\sigma$-compact, the positive measures in both parts are
Radon measures in the sense of
Definition~\ref{def:b4-espacios-medida-localmente-compacto}.
This applies, in particular, to open subsets of $\mathbb R^n$ and to
second-countable smooth manifolds.
\end{theorem}

\begin{proof}
The positive representation in part (a), with its regularity and
uniqueness properties, is Riesz's theorem; see
\cite{Bogachev2007,Cohn2013}. We shall deduce part (b) from it without assuming
that one of the representing measures has finite total mass.

For $f\in C_c(X)$ with $f\geq0$, set
\[
L^+(f):=\sup\{L(g)\mid g\in C_c(X),\ 0\leq g\leq f\}.
\]
If $K=\supp(f)$, every admissible function $g$ is supported in $K$ and
$\|g\|_\infty\leq\|f\|_\infty$. By local boundedness of $L$,
\[
0\leq L^+(f)\leq c_K\|f\|_\infty<\infty.
\]
Homogeneity for nonnegative scalars follows directly from
the definition, including the zero scalar.

Let $f,h\in C_c(X)$ be nonnegative. If $0\leq g_1\leq f$ and
$0\leq g_2\leq h$, then $0\leq g_1+g_2\leq f+h$. Taking suprema
gives $L^+(f+h)\geq L^+(f)+L^+(h)$. For the reverse inequality,
let $0\leq g\leq f+h$ and define
\[
g_1:=\min\{g,f\},\qquad g_2:=g-g_1.
\]
These functions are continuous and compactly supported, and satisfy
$0\leq g_1\leq f$ and $0\leq g_2\leq h$. Thus
$L(g)=L(g_1)+L(g_2)\leq L^+(f)+L^+(h)$.
Taking the supremum over $g$ shows that $L^+$ is additive on the positive
cone.

Extend $L^+$ to all of $C_c(X)$ by
$L^+(f-h):=L^+(f)-L^+(h)$ for $f,h\geq0$.
If $f-h=f_1-h_1$, then $f+h_1=f_1+h$, and the additivity just proved
shows that both expressions give the same value. We thus obtain a
positive linear functional. Since $g=f$ is admissible in the supremum when
$f\geq0$, we have $L^+(f)\geq L(f)$. Consequently,
$L^-:=L^+-L$ is also a positive linear functional.
Applying (a) to $L^+$ and $L^-$ yields the measures $\mu_1,\mu_2$
with the asserted representation. Both integrals are finite for each
$u\in C_c(X)$ because the measures are finite on its support.

Now let $U\subseteq X$ be open with compact closure.
Then $\mu_i(U)\leq\mu_i(\overline U)<\infty$ for $i=1,2$.
The restrictions of $\mu_i$ to $U$ are regular: all their Borel sets
have finite measure, and outer regularity is preserved by intersecting
the approximating open sets with $U$. The difference
\[
\nu_U:=\mu_1\restriction_U-\mu_2\restriction_U
\]
is a finite signed measure. To check that its variations are
regular, set $\eta:=\mu_1\restriction_U+\mu_2\restriction_U$.
The Hahn decomposition gives $0\leq\nu_U^\pm\leq\eta$.
Given a Borel set $E\subseteq U$ and $\varepsilon>0$, regularity and
finiteness of $\eta$ allow us to choose a compact set $K\subseteq E$ and an open set
$O\supseteq E$ in $U$ such that $\eta(O\setminus K)<\varepsilon$.
Then
\[
0\leq\nu_U^\pm(O)-\nu_U^\pm(K)
=\nu_U^\pm(O\setminus K)<\varepsilon.
\]
This proves regularity of both variations. Thus $\nu_U$ is a
finite signed Radon measure.
These differences are compatible under restriction from one open set to another,
so they define a local signed Radon measure $\nu$.

For uniqueness, let $\nu$ and $\widetilde\nu$ be two representing families,
and fix $U$ as above. A function in $C_c(U)$ extends by zero to
an element of $C_c(X)$, so the finite measure
$\lambda:=\nu_U-\widetilde\nu_U$ vanishes on $C_c(U)$.
By Hahn--Jordan,
$\int_Uu\,d\lambda^+=\int_Uu\,d\lambda^-$ for every $u\in C_c(U)$.
Uniqueness in (a) implies $\lambda^+=\lambda^-$; since they are mutually
singular, both vanish. Thus $\nu_U=\widetilde\nu_U$ for every $U$.

Finally, suppose that $X$ is $\sigma$-compact, and let $\mu$ be any
of the positive measures obtained. Choose increasing compact sets
$K_j$ covering $X$. For a Borel set $E$, continuity from below gives
\[
\mu(E)=\lim_{j\to\infty}\mu(E\cap K_j).
\]
Each $E\cap K_j$ has finite measure and is inner regular by (a).
Hence $\mu(E)=\displaystyle\sup\{\mu(K)\mid K\subseteq E,\ K\text{ compact}\}$,
also when $\mu(E)=\infty$. Outer regularity is already part
of (a), and finiteness on compact sets implies local finiteness on a locally
compact space. This proves the last assertion.
\end{proof}

The local family $\nu$ must not be confused with a difference defined
on all Borel subsets of $X$: when
$\mu_1(E)=\mu_2(E)=\infty$, the expression $\mu_1(E)-\mu_2(E)$ is not
meaningful. For example, $L(u)=\int_{\mathbb R}xu(x)\,dx$ is locally
bounded on $C_c(\mathbb R)$ and is represented by
$\displaystyle\max\{x,0\}\,dx$ and $\displaystyle\max\{-x,0\}\,dx$, both of infinite total mass.
On each relatively compact open set, however, both measures are
finite and their difference is well defined.

\section{Lebesgue and Bochner integration}\label{ap:integracion-de-lebesgue-y-bochner}

Having fixed the basic notions of measure, we recall the construction of the Lebesgue integral and its extension to Banach-valued functions through the Bochner integral. We begin with the vector-valued Riemann integral, which provides a natural way to introduce compatibility of integration with linear operators. For a more detailed treatment of the Bochner integral, see \cite{DrabekMilota2013}.

\subsection{The Banach-valued Riemann integral}

\begin{definition}\label{def:b5-integral-riemann-banach}\index{Riemann integral!Banach-valued}
Let $E$ be a Banach space, let $a,b\in\mathbb{R}$ with $a<b$, and let $f\colon [a,b]\longrightarrow E$. A \textbf{tagged partition} of $[a,b]$ is a pair consisting of a partition $\mathcal{P}=\{a=t_{0}<t_{1}<\cdots<t_{n}=b\}$ and points $\tau_{i}\in[t_{i-1},t_{i}]$, $i\in\{1,\dots,n\}$. Define the mesh of the partition by $\displaystyle|\mathcal{P}|:=\displaystyle\max_{1\leq i\leq n}(t_{i}-t_{i-1})$. We say that $f$ is \textbf{Riemann integrable} if there exists $I\in E$ such that, for each $\varepsilon>0$, there exists $\delta>0$ with the property that $\displaystyle\left\|\displaystyle\sum_{i=1}^{n}f(\tau_{i})(t_{i}-t_{i-1})-I\right\|_{E}<\varepsilon$ for every tagged partition with $|\mathcal{P}|<\delta$. In this case, $I$ is unique and is denoted by $\displaystyle\int_{a}^{b}f(t)\,dt$.
\end{definition}

\begin{theorem}[Graves]\label{teo:b5-graves-integral-riemann-banach}
Let $E$ be a Banach space. Every continuous function $f\colon [a,b]\longrightarrow E$ is Riemann integrable.
\end{theorem}
\begin{proof}
Since $[a,b]$ is compact, $f$ is uniformly continuous. For a tagged partition \((\mathcal P,\tau)\), denote its Riemann sum by
\[
S(f;\mathcal P,\tau)
:=
\sum_{i=1}^{n}f(\tau_i)(t_i-t_{i-1})
\]
Fix \(\varepsilon>0\). Uniform continuity allows us to choose \(\delta>0\) so that
\[
|s-t|<\delta
\quad\Longrightarrow\quad
\|f(s)-f(t)\|_E<\frac{\varepsilon}{4(b-a)}.
\]
Let \((\mathcal P,\tau)\) and \((\mathcal Q,\sigma)\) be two tagged partitions with mesh less than \(\delta\), and let \(\mathcal R\) be the partition obtained by combining and ordering the points of \(\mathcal P\) and \(\mathcal Q\). Choose a tag in each subinterval of \(\mathcal R\). Each such subinterval is contained in an interval of \(\mathcal P\); its tag and the corresponding tag of \(\mathcal P\) are less than \(\delta\) apart. By the triangle inequality,
\[
\bigl\|S(f;\mathcal P,\tau)-S(f;\mathcal R,\rho)\bigr\|_E
\leq
\frac{\varepsilon}{4(b-a)}
\sum_{J\in\mathcal R}|J|
=\frac{\varepsilon}{4}.
\]
The same estimate holds with \(\mathcal Q\) in place of \(\mathcal P\). Consequently,
\[
\bigl\|S(f;\mathcal P,\tau)-S(f;\mathcal Q,\sigma)\bigr\|_E
<\frac{\varepsilon}{2}.
\]

For each \(k\in\mathbb N\), choose a tagged partition \((\mathcal P_k,\tau^{(k)})\) with mesh less than \(\displaystyle\frac{1}{k}\), and set \(S_k:=S(f;\mathcal P_k,\tau^{(k)})\). The preceding estimate shows that \((S_k)\) is Cauchy. Since \(E\) is complete, there exists \(I\in E\) such that \(S_k\to I\).

It remains to check that all sufficiently fine sums converge to \(I\). Given \(\varepsilon>0\), repeat the preceding argument with \(\varepsilon\) in place of the quantity fixed at the beginning and obtain the corresponding \(\delta\). Choose \(k\) so large that \(|\mathcal P_k|<\delta\) and \(\|S_k-I\|_E<\displaystyle\frac{\varepsilon}{2}\). For any tagged partition \((\mathcal P,\tau)\) with \(|\mathcal P|<\delta\), comparison through a common refinement gives
\[
\|S(f;\mathcal P,\tau)-I\|_E
\leq
\|S(f;\mathcal P,\tau)-S_k\|_E+\|S_k-I\|_E
<\varepsilon.
\]
This is precisely Riemann integrability of \(f\).
\end{proof}

\begin{proposition}\label{prop:b5-propiedades-integral-riemann-banach}
Let $E$ be a Banach space, and let $f,g\colon [a,b]\longrightarrow E$ be Riemann integrable functions. Then $\alpha f+\beta g$ is Riemann integrable for any $\alpha,\beta\in\mathbb{K}$ and $\displaystyle\int_{a}^{b}(\alpha f+\beta g)(t)\,dt=\alpha\int_{a}^{b}f(t)\,dt+\beta\int_{a}^{b}g(t)\,dt$. If, in addition, the scalar function $t\longmapsto\|f(t)\|_{E}$ is Riemann integrable, then $\displaystyle\left\|\int_{a}^{b}f(t)\,dt\right\|_{E}\leq\int_{a}^{b}\|f(t)\|_{E}\,dt$. In particular, the map $\displaystyle\mathcal{I}\colon C([a,b],E)\longrightarrow E$, given by $\displaystyle\mathcal{I}(f):=\int_{a}^{b}f(t)\,dt$, is linear and bounded, and satisfies $\|\mathcal{I}\|_{\mathcal L(C([a,b],E),E)}\leq b-a$ with respect to the norm $\displaystyle\|f\|_{C([a,b],E)}:=\displaystyle\max_{t\in[a,b]}\|f(t)\|_{E}$.
\end{proposition}
\begin{proof}
Linearity follows by taking linear combinations of Riemann sums. For the inequality, suppose that the scalar function $t\longmapsto\|f(t)\|_E$ is Riemann integrable. Every Riemann sum satisfies $\displaystyle\left\|\displaystyle\sum_{i=1}^{n}f(\tau_{i})(t_{i}-t_{i-1})\right\|_{E}\leq\displaystyle\sum_{i=1}^{n}\|f(\tau_{i})\|_{E}(t_{i}-t_{i-1})$; passing to the limit, using integrability of both functions, gives the assertion. If $f\in C([a,b],E)$, both $f$ and $t\longmapsto\|f(t)\|_E$ are continuous and hence Riemann integrable by Theorem~\ref{teo:b5-graves-integral-riemann-banach}. Finally, $\displaystyle\int_{a}^{b}\|f(t)\|_{E}\,dt\leq(b-a)\|f\|_{C([a,b],E)}$.
\end{proof}

The integral also commutes with continuous linear operators and, under an additional integrability hypothesis, with closed linear operators.

\begin{proposition}\label{prop:b5-integral-riemann-operadores-lineales}
Let $E$ and $F$ be Banach spaces, and let $f\colon [a,b]\longrightarrow E$ be Riemann integrable.
\begin{enumerate}[label=(\alph*)]
\item If $T\in\mathcal{L}(E,F)$, then $T\circ f$ is Riemann integrable and $\displaystyle T\left(\int_{a}^{b}f(t)\,dt\right)=\int_{a}^{b}T(f(t))\,dt$.
\item If $T\colon \operatorname{Dom}(T)\subseteq E\longrightarrow F$ is a closed linear operator, $f(t)\in\operatorname{Dom}(T)$ for every $t\in[a,b]$, and $T\circ f$ is Riemann integrable, then $\displaystyle\int_{a}^{b}f(t)\,dt\in\operatorname{Dom}(T)$ and the same equality holds.
\end{enumerate}
\end{proposition}
\begin{proof}
The first assertion follows by applying $T$ to Riemann sums and using continuity. For the second, let $(S_{k})_{k\in\mathbb{N}}$ be a sequence of Riemann sums of $f$ associated with partitions whose mesh tends to zero. Since $\operatorname{Dom}(T)$ is a vector subspace, $S_{k}\in\operatorname{Dom}(T)$ and $T(S_{k})$ is the corresponding Riemann sum for $T\circ f$. Thus $\displaystyle S_{k}\longrightarrow\int_{a}^{b}f(t)\,dt$ in $E$ and $\displaystyle T(S_{k})\longrightarrow\int_{a}^{b}T(f(t))\,dt$ in $F$. Closedness of $T$ completes the proof.
\end{proof}

\begin{theorem}[Fundamental theorem of calculus for Banach-valued functions]
\label{teo:b5-fundamental-calculo-riemann-banach}
\index{fundamental theorem of calculus!for Banach-valued functions}
Let $E$ be a Banach space, and let $F\colon [a,b]\longrightarrow E$ be a function of class $C^1$. Then
\[
F(b)-F(a)=\int_a^bF'(t)\,dt.
\]
More generally, for any $s,t\in[a,b]$ with $s\leq t$, we have $\displaystyle F(t)-F(s)=\int_s^tF'(\tau)\,d\tau$.
\end{theorem}

\begin{proof}
Continuity of $F'$ and Theorem~\ref{teo:b5-graves-integral-riemann-banach} ensure that the integral is well defined. Let $\varphi\in E'$. The scalar function $\varphi\circ F$ is of class $C^1$ and, by the chain rule, $(\varphi\circ F)'(t)=\varphi(F'(t))$. The fundamental theorem of calculus for real-valued functions and Proposition~\ref{prop:b5-integral-riemann-operadores-lineales} imply
\[
\varphi(F(b)-F(a))
=
\int_a^b\varphi(F'(t))\,dt
=
\varphi\left(\int_a^bF'(t)\,dt\right).
\]
Thus $\displaystyle\varphi\left(F(b)-F(a)-\int_a^bF'(t)\,dt\right)=0$ for every $\varphi\in E'$. Since the continuous dual separates points, the identity follows. The formula on $[s,t]$ is proved in the same way.
\end{proof}

\begin{proposition}[Integration by parts for Banach-valued functions]
\label{prop:b5-integracion-por-partes-riemann-banach}
\index{integration by parts!for Banach-valued functions}
Let $E$ be a Banach space, let $\alpha\colon [a,b]\longrightarrow\mathbb K$ be of class $C^1$, and let $F\colon [a,b]\longrightarrow E$ be of class $C^1$. Then
\begin{equation}
\label{eq:b5-integracion-por-partes-riemann-banach}
\int_a^b\alpha(t)F'(t)\,dt
=
\alpha(b)F(b)-\alpha(a)F(a)-\int_a^b\alpha'(t)F(t)\,dt.
\end{equation}
\end{proposition}

\begin{proof}
The map $B\colon \mathbb K\times E\longrightarrow E$, defined by $B(\lambda,v):=\lambda v$, is continuous and bilinear. The product rule gives
\[
\frac{d}{dt}\bigl(\alpha(t)F(t)\bigr)
=
\alpha'(t)F(t)+\alpha(t)F'(t).
\]
Applying Theorem~\ref{teo:b5-fundamental-calculo-riemann-banach} to the function $t\longmapsto\alpha(t)F(t)$ and using linearity of the integral, we obtain
\[
\alpha(b)F(b)-\alpha(a)F(a)
=
\int_a^b\alpha'(t)F(t)\,dt+
\int_a^b\alpha(t)F'(t)\,dt.
\]
Solving for the last integral gives \eqref{eq:b5-integracion-por-partes-riemann-banach}.
\end{proof}

\subsection{The Lebesgue integral}

From now on, $(M,\mathcal{A},\mu)$ denotes a measure space. Measurability of a function will always be understood with respect to the $\sigma$-algebra $\mathcal{A}$ on the domain and the Borel $\sigma$-algebra on the codomain.

\begin{definition}\label{def:b5-integral-lebesgue-medida}\index{measurable function}
Let $Y$ be a topological space. A function $f\colon M\longrightarrow Y$ is \textbf{$\mathcal{A}$-measurable}, or simply \textbf{measurable}, if $f^{-1}(U)\in\mathcal{A}$ for every open set $U\subseteq Y$. In particular, for functions $f\colon M\longrightarrow\overline{\mathbb{R}}:=[-\infty,\infty]$, we equip $\overline{\mathbb{R}}$ with the order topology.
\end{definition}

The Lebesgue integral is first constructed for nonnegative simple functions. We use a single notion of simple function, valid for both scalar-valued and vector-valued functions.

\begin{definition}\label{def:b5-integral-lebesgue-espacio-medible-ajenos-dos-dos-funcion}\index{simple function}
Let $V$ be a vector space. A function $s\colon M\longrightarrow V$ is \textbf{simple} if there exist pairwise disjoint measurable sets $A_{1},\dots,A_{r}\in\mathcal{A}$ and vectors $v_{1},\dots,v_{r}\in V$ such that $\displaystyle s=\displaystyle\sum_{i=1}^{r}v_{i}\mathbf{1}_{A_{i}}$, where $\mathbf{1}_{A}$ denotes the indicator function of $A$.
\end{definition}

\begin{definition}\label{def:b5-integral-lebesgue-integral-espacio-medida-funcion-simple}\index{integral of a simple function}
Let $s\colon M\longrightarrow[0,\infty)$ be a simple function, and write $\displaystyle s=\displaystyle\sum_{i=1}^{r}a_{i}\mathbf{1}_{A_{i}}$, where $a_{i}\geq0$ and the sets $A_{i}$ are measurable and pairwise disjoint. Define $\displaystyle\int_{M}s\,d\mu:=\displaystyle\sum_{i=1}^{r}a_{i}\mu(A_{i})$, with the convention $0\cdot\infty=0$. This definition is independent of the chosen representation of $s$.
\end{definition}

\begin{definition}\label{def:b5-integral-lebesgue-integral}\index{Lebesgue integral!of a nonnegative measurable function}
Let $f\colon M\longrightarrow[0,\infty]$ be measurable. Define its \textbf{Lebesgue integral} by $\displaystyle\int_{M}f\,d\mu:=\displaystyle\sup\left\{\int_{M}s\,d\mu\ \middle|\ s\colon M\longrightarrow[0,\infty)\text{ is simple and }0\leq s\leq f\right\}$. If $A\in\mathcal{A}$, define $\displaystyle\int_{A}f\,d\mu:=\int_{M}f\mathbf{1}_{A}\,d\mu$.
\end{definition}

The first fundamental convergence theorem expresses continuity of the integral with respect to increasing limits.

\begin{theorem}[Monotone convergence theorem]\label{convergencia monotona}\index{convergence theorem!monotone}
Let $f\colon M\longrightarrow[0,\infty]$ be measurable, and let $(f_{n})_{n\in\mathbb{N}}$ be a sequence of measurable functions $f_{n}\colon M\longrightarrow[0,\infty]$ such that $f_{n}\leq f_{n+1}$ for every $n\in\mathbb{N}$ and $f_{n}(x)\longrightarrow f(x)$ for $\mu$-almost every $x\in M$. Then $\displaystyle\int_{A}f\,d\mu=\lim_{n\to\infty}\int_{A}f_{n}\,d\mu$ for every $A\in\mathcal{A}$. The measurability assumption on $f$ is automatic if $\mu$ is complete.
\end{theorem}
\begin{proof}
The function $h:=\displaystyle\sup_{n\in\mathbb N}f_n$ is measurable and agrees with $f$ outside a measurable null set $N$. If $\mu$ is complete, then for every open set $U\subseteq[0,\infty]$ we have $f^{-1}(U)=(h^{-1}(U)\setminus N)\cup(f^{-1}(U)\cap N)$, and both sets in this union are measurable. This proves the last assertion of the statement. In the general case, measurability of $f$ is assumed. Replacing $f$ and all the $f_n$ by zero on $N$ preserves their measurability and integrals and gives convergence and monotonicity at every point. It suffices to prove the equality for \(A=M\), since the general case follows by applying the same argument to \(f_n\mathbf 1_A\).

By monotonicity of the integral, the sequence \(\displaystyle \int_Mf_n\,d\mu\) is increasing and
\[
I:=\lim_{n\to\infty}\int_Mf_n\,d\mu
\leq\int_Mf\,d\mu.
\]
For the reverse inequality, let \(s\) be a nonnegative simple function with \(s\leq f\), and fix \(c\in(0,1)\). Define
\[
A_n:=\{x\in M\mid f_n(x)\geq c\,s(x)\}.
\]
The sets \(A_n\) are measurable, form an increasing sequence, and have union \(M\): if \(s(x)=0\), then \(x\in A_n\) for every \(n\); if \(s(x)>0\), convergence \(f_n(x)\to f(x)\geq s(x)>c\,s(x)\) implies that \(x\in A_n\) for sufficiently large \(n\). Thus
\[
\int_Mf_n\,d\mu
\geq c\int_{A_n}s\,d\mu.
\]
If \(\displaystyle s=\displaystyle\sum_{j=1}^r a_j\mathbf1_{B_j}\) with \(a_j\geq0\) and the \(B_j\) pairwise disjoint, continuity from below of the measure gives
\[
\lim_{n\to\infty}\int_{A_n}s\,d\mu
=
\sum_{j=1}^r a_j\lim_{n\to\infty}\mu(A_n\cap B_j)
=
\sum_{j=1}^r a_j\mu(B_j)
=\int_Ms\,d\mu.
\]
Consequently, \(\displaystyle I\geq c\int_Ms\,d\mu\). Letting \(c\to1^{-}\) and then taking the supremum over all simple functions \(0\leq s\leq f\) gives
\[
I\geq\int_Mf\,d\mu.
\]
Together with the first inequality, this proves equality.
\end{proof}

\begin{lemma}[Fatou]\label{lem:fatou}
\index{Fatou's lemma}
Let $(f_n)_{n\in\mathbb N}$ be a sequence of measurable functions $f_n\colon M\longrightarrow[0,\infty]$. Then
\[
\int_M\liminf_{n\to\infty}f_n\,d\mu
\leq
\liminf_{n\to\infty}\int_M f_n\,d\mu.
\]
\end{lemma}

\begin{proof}
For each $n\in\mathbb N$, define $g_n:=\displaystyle\inf_{k\geq n}f_k$. The sequence $(g_n)$ is increasing and converges pointwise to $\displaystyle\liminf_{n\to\infty}f_n$. By Theorem~\ref{convergencia monotona},
\[
\int_M\liminf_{n\to\infty}f_n\,d\mu
=
\lim_{n\to\infty}\int_M g_n\,d\mu.
\]
Since $g_n\leq f_k$ for every $k\geq n$, we have $\displaystyle\int_Mg_n\,d\mu\leq\inf_{k\geq n}\int_Mf_k\,d\mu$. Letting $n\to\infty$ yields the desired inequality.
\end{proof}

\begin{definition}\label{def:b5-integral-lebesgue-lebesgue-integrable}\index{Lebesgue integrable}
Let $f\colon M\longrightarrow\overline{\mathbb{R}}$ be measurable, and define $f^{+}:=\displaystyle\max\{f,0\}$ and $f^{-}:=\displaystyle\max\{-f,0\}$. Whenever at least one of the quantities $\displaystyle\int_{M}f^{+}\,d\mu$ and $\displaystyle\int_{M}f^{-}\,d\mu$ is finite, define $\displaystyle\int_{M}f\,d\mu:=\int_{M}f^{+}\,d\mu-\int_{M}f^{-}\,d\mu$. We say that $f$ is \textbf{Lebesgue integrable}, or simply \textbf{integrable}, if $\displaystyle\int_{M}|f|\,d\mu<\infty$. Denote by $\mathcal{L}^{1}(M,\mu)$ the vector space of integrable functions and by $L^{1}(M,\mu)$ the space of their equivalence classes modulo equality almost everywhere. For $A\in\mathcal{A}$, write $\displaystyle\int_{A}f\,d\mu:=\int_{M}f\mathbf{1}_{A}\,d\mu$.
\end{definition}

\begin{theorem}[Lebesgue dominated convergence theorem]\label{convergencia dominada}\index{convergence theorem!Lebesgue dominated}
Let $f\colon M\longrightarrow\overline{\mathbb{R}}$ be measurable, and let $(f_{n})_{n\in\mathbb{N}}$ be a sequence of measurable functions $f_{n}\colon M\longrightarrow\overline{\mathbb{R}}$ such that $f_{n}(x)\longrightarrow f(x)$ for $\mu$-almost every $x\in M$. Suppose that there exists $g\in\mathcal{L}^{1}(M,\mu)$ such that $|f_{n}|\leq g$ $\mu$-almost everywhere for every $n\in\mathbb{N}$. Then $f\in\mathcal{L}^{1}(M,\mu)$, $\displaystyle\int_{M}|f_{n}-f|\,d\mu\longrightarrow0$, and $\displaystyle\int_{A}f\,d\mu=\lim_{n\to\infty}\int_{A}f_{n}\,d\mu$ for every $A\in\mathcal{A}$.
\end{theorem}
\begin{proof}
The inequality \(|f_n|\leq g\) implies \(g\geq0\) almost everywhere;
if desired, \(g\) may be replaced by \(|g|\). We may modify the
functions on a common null set and assume that convergence and
domination hold at every point. The pointwise limit \(f\) is measurable and
satisfies \(|f|\leq g\), so \(f\in\mathcal L^1(M,\mu)\). Moreover,
\[
0\leq |f_n-f|\leq2g
\quad\text{and}\quad
2g-|f_n-f|\longrightarrow2g
\]
pointwise. Applying Fatou's lemma to the nonnegative functions \(2g-|f_n-f|\) gives
\[
2\int_Mg\,d\mu
\leq
\liminf_{n\to\infty}
\left(2\int_Mg\,d\mu-\int_M|f_n-f|\,d\mu\right)
=
2\int_Mg\,d\mu-
\limsup_{n\to\infty}\int_M|f_n-f|\,d\mu.
\]
Since \(\displaystyle \int_Mg\,d\mu<\infty\), it follows that
\[
\lim_{n\to\infty}\int_M|f_n-f|\,d\mu=0.
\]
Finally, for every \(A\in\mathcal A\),
\[
\left|\int_Af_n\,d\mu-\int_Af\,d\mu\right|
\leq\int_A|f_n-f|\,d\mu
\leq\int_M|f_n-f|\,d\mu\longrightarrow0.
\]
\end{proof}

In Euclidean spaces, the Lebesgue integral transforms through the Jacobian of a diffeomorphism.

\begin{theorem}[Change-of-variables theorem]\label{teorema de cambio de variable}\index{change-of-variables theorem}
Let $\Omega,\Omega'\subseteq\mathbb{R}^{n}$ be open sets, and let $\phi\colon \Omega\longrightarrow\Omega'$ be a diffeomorphism of class $C^{1}$.
\begin{enumerate}[label=(\alph*)]
\item If $f\colon \Omega'\longrightarrow[0,\infty]$ is measurable, then $\displaystyle\int_{\Omega}f(\phi(x))|\det D\phi(x)|\,d\lambda_{n}(x)=\int_{\Omega'}f(y)\,d\lambda_{n}(y)$.
\item A measurable function $f\colon \Omega'\longrightarrow\mathbb{R}$ is integrable if and only if $(f\circ\phi)|\det D\phi|$ is integrable on $\Omega$, and in this case the same equality holds.
\end{enumerate}
\end{theorem}

A proof based on approximation by cubes and the local behavior of a \(C^1\) diffeomorphism can be found in \cite{Grabinsky,Bogachev2007}.

\begin{corollary}\label{corolario del cambio de variable}
Let $\Omega\subseteq\mathbb{R}^{n}$ be open, and let $A\in\mathcal{L}(\mathbb{R}^{n})$ satisfy $A\subseteq\Omega$ and $\lambda_{n}(A\setminus\operatorname{int}(A))=0$. Suppose that $\Phi\in C^{1}(\Omega,\mathbb{R}^{n})$ and that $\Phi\restriction_{\operatorname{int}(A)}\colon \operatorname{int}(A)\longrightarrow\Phi(\operatorname{int}(A))$ is a diffeomorphism.
\begin{enumerate}[label=(\alph*)]
\item If $f\colon \Phi(A)\longrightarrow[0,\infty]$ is measurable, then $\displaystyle\int_{\Phi(A)}f(y)\,d\lambda_{n}(y)=\int_{A}f(\Phi(x))|\det D\Phi(x)|\,d\lambda_{n}(x)$.
\item A measurable function $f\colon \Phi(A)\longrightarrow\mathbb{R}$ is integrable if and only if $(f\circ\Phi)|\det D\Phi|$ is integrable on $A$, and in this case the same equality holds.
\end{enumerate}
\end{corollary}
\begin{proof}
Let $N:=A\setminus\operatorname{int}(A)$. By Proposition~\ref{imagen de lebesgue medible bajo c1}, $\Phi(N)$ is Lebesgue measurable and $\lambda_{n}(\Phi(N))=0$, since $\Phi$ is locally Lipschitz. Moreover, $\Phi(A)\setminus\Phi(\operatorname{int}(A))\subseteq\Phi(N)$. Thus the integrals over $A$ and $\Phi(A)$ do not change when these sets are replaced by $\operatorname{int}(A)$ and $\Phi(\operatorname{int}(A))$, respectively. The result follows from Theorem~\ref{teorema de cambio de variable}.
\end{proof}

\begin{theorem}[Tonelli's theorem]\label{teo:tonelli}\index{Tonelli's theorem}
Let $(X,\mathcal A,\mu)$ and $(Y,\mathcal B,\nu)$ be $\sigma$-finite measure spaces. If $f\colon X\times Y\longrightarrow[0,\infty]$ is $\mathcal A\otimes\mathcal B$-measurable, then the functions
\[
x\longmapsto\int_Y f(x,y)\,d\nu(y),
\qquad
y\longmapsto\int_X f(x,y)\,d\mu(x)
\]
are measurable and
\[
\int_{X\times Y}f\,d(\mu\otimes\nu)
=
\int_X\left(\int_Yf(x,y)\,d\nu(y)\right)d\mu(x)
=
\int_Y\left(\int_Xf(x,y)\,d\mu(x)\right)d\nu(y),
\]
with values in $[0,\infty]$ allowed. The notation $\mu\otimes\nu$
denotes the unique measure on $\mathcal A\otimes\mathcal B$ with
$(\mu\otimes\nu)(A\times B)=\mu(A)\nu(B)$ for measurable rectangles,
with $0\cdot\infty=0$; its existence is included in the proof.
The completed version has the same identities, with sections and partial
integrals defined almost everywhere in the completed factor spaces.
\end{theorem}
\begin{proof}
We first recall the elementary $\pi$--$\lambda$ principle. A $\lambda$-system
contains the whole space and is closed under complements and countable
disjoint unions; a $\pi$-system is closed under finite intersections.
The smallest $\lambda$-system $\mathcal D$ containing a $\pi$-system
$\mathcal P$ is its generated $\sigma$-algebra. Indeed, for fixed
$A\in\mathcal P$, the sets $B\in\mathcal D$ for which $A\cap B\in\mathcal D$
form a $\lambda$-system containing $\mathcal P$, hence all of $\mathcal D$.
Repeating with $A\in\mathcal D$ shows closure under finite intersections.
A $\lambda$-system with that property is a $\sigma$-algebra, by disjointizing
countable unions. This proves the principle.

For a product-measurable set $E$, put
$E_x:=\{y\in Y\mid(x,y)\in E\}$ and
$E^y:=\{x\in X\mid(x,y)\in E\}$.
All sections are measurable: sets having measurable sections form a
$\sigma$-algebra containing the rectangles. Suppose first that both
measures are finite. Let $\mathcal D$ be the class of product-measurable
sets $E$ such that $x\mapsto\nu(E_x)$ and $y\mapsto\mu(E^y)$ are measurable
and their integrals are equal. Rectangles belong to $\mathcal D$.
Complements belong to it because the partial measures subtract from the
finite constants $\nu(Y)$ and $\mu(X)$, and the total integrals subtract
from $\mu(X)\nu(Y)$. For a disjoint sequence, countable additivity on
sections and monotone convergence prove the same assertions for its union.
Thus $\mathcal D$ is a $\lambda$-system. The $\pi$--$\lambda$ principle
gives $\mathcal D=\mathcal A\otimes\mathcal B$.

Define
\[
 \lambda(E):=\int_X\nu(E_x)\,d\mu(x)
            =\int_Y\mu(E^y)\,d\nu(y).
\]
Monotone convergence for disjoint unions shows that $\lambda$ is a measure,
and its values on rectangles are the required products. Two such finite
measures agree on the generated $\sigma$-algebra by the same
$\pi$--$\lambda$ argument, so the measure is unique.

For $\sigma$-finite spaces, choose increasing sequences $(X_j)_{j\in\mathbb N}$ and $(Y_j)_{j\in\mathbb N}$ of measurable sets
of finite measure, with $\displaystyle X=\bigcup_{j=1}^{\infty}X_j$ and
$\displaystyle Y=\bigcup_{j=1}^{\infty}Y_j$. Apply the finite case
on each $X_j\times Y_j$. Its measures are consistent on smaller blocks
by uniqueness. Consequently
\[
 (\mu\otimes\nu)(E)
 :=\lim_{j\to\infty}
       \int_{X_j}\nu(E_x\cap Y_j)\,d\mu(x)
\]
defines a countably additive measure; interchange of the increasing limit
with a sum of nonnegative terms proves countable additivity.
It is unique because the blocks exhaust the product. Monotone convergence
on these blocks proves measurability of the partial measures and both
iterated-integral identities for $1_E$. This is where $\sigma$-finiteness
enters the argument.

Linearity now proves the identities for arbitrary nonnegative simple
functions, whose level sets need not be rectangles. Approximate $f$
monotonically by such functions and apply
Theorem~\ref{convergencia monotona} to the inner, outer, and product
integrals. Their increasing limits prove the asserted equalities and
measurability for $f$.

Finally, a function measurable for the completed product measure has a
product-measurable representative $f_0$ agreeing with it outside a subset
of a product-measurable null set $N$. The indicator case just proved gives
$\nu(N_x)=0$ for almost every $x$ and $\mu(N^y)=0$ for almost every $y$.
In the completed factor spaces the corresponding sections of $f$ and
$f_0$ are therefore measurable and have equal integrals. This proves the
completed version, including its almost-everywhere qualifications.
\end{proof}

\begin{theorem}[Fubini's theorem]\label{Fubini}\index{Fubini's theorem}
Let $1\leq m<n$, and identify $\mathbb{R}^{n}=\mathbb{R}^{m}\times\mathbb{R}^{n-m}$. If $f\in\mathcal{L}^{1}(\mathbb{R}^{n},\lambda_{n})$, then, for $\lambda_{n-m}$-almost every $y\in\mathbb{R}^{n-m}$, the function $x\longmapsto f(x,y)$ belongs to $\mathcal{L}^{1}(\mathbb{R}^{m},\lambda_{m})$. The function $\displaystyle F(y):=\int_{\mathbb{R}^{m}}f(x,y)\,d\lambda_{m}(x)$, defined arbitrarily on the exceptional set, belongs to $\mathcal{L}^{1}(\mathbb{R}^{n-m},\lambda_{n-m})$, and $\displaystyle\int_{\mathbb{R}^{n}}f\,d\lambda_{n}=\int_{\mathbb{R}^{n-m}}F(y)\,d\lambda_{n-m}(y)$. Similarly, for $\lambda_{m}$-almost every $x\in\mathbb{R}^{m}$, the function $y\longmapsto f(x,y)$ is integrable. If $\displaystyle H(x):=\int_{\mathbb{R}^{n-m}}f(x,y)\,d\lambda_{n-m}(y)$, then $H\in\mathcal{L}^{1}(\mathbb{R}^{m},\lambda_{m})$ and $\displaystyle\int_{\mathbb{R}^{n}}f\,d\lambda_{n}=\int_{\mathbb{R}^{m}}H(x)\,d\lambda_{m}(x)$.
\end{theorem}
\begin{proof}
The completion of $\lambda_m\otimes\lambda_{n-m}$ is $\lambda_n$:
both agree on Borel rectangles, uniqueness identifies their Borel
restrictions, and completion adds exactly the subsets of null sets.
Thus the completed version of Tonelli's theorem applies to $|f|$ and
ensures measurability of almost every section. Since
\[
\int_{\mathbb R^{n-m}}
\left(\int_{\mathbb R^m}|f(x,y)|\,d\lambda_m(x)\right)
d\lambda_{n-m}(y)
=\int_{\mathbb R^n}|f|\,d\lambda_n<\infty,
\]
the inner integral of \(|f|\) is finite for \(\lambda_{n-m}\)-almost every \(y\). For these \(y\), the function \(x\mapsto f(x,y)\) is integrable. Writing \(f=f^+-f^-\) and applying Tonelli separately to \(f^+\) and \(f^-\) gives measurability of
\[
F(y)=\int_{\mathbb R^m}f(x,y)\,d\lambda_m(x)
\]
outside the exceptional set and the identity
\[
\int_{\mathbb R^{n-m}}F(y)\,d\lambda_{n-m}(y)
=
\int_{\mathbb R^n}f^+\,d\lambda_n
-
\int_{\mathbb R^n}f^-\,d\lambda_n
=
\int_{\mathbb R^n}f\,d\lambda_n.
\]
Likewise,
\[
|F(y)|
\leq\int_{\mathbb R^m}|f(x,y)|\,d\lambda_m(x)
\]
for almost every \(y\), and the right-hand side is integrable; thus \(F\in\mathcal L^1(\mathbb R^{n-m})\). Interchanging the roles of \(x\) and \(y\) in this argument gives the assertions about \(H\). For complex-valued $f$, apply the real-valued argument to its real and imaginary parts.
\end{proof}

Before stating the consequence ensuring uniqueness of weak derivatives, we define local integrability.

\begin{definition}\label{def:b5-integral-lebesgue-integrable-localmente}\index{locally integrable function}
Let $\Omega\subseteq\mathbb{R}^{n}$ be open. A measurable function $f\colon \Omega\longrightarrow\overline{\mathbb{R}}$ is \textbf{locally integrable} if $\displaystyle\int_{K}|f|\,d\lambda_{n}<\infty$ for every compact set $K\subseteq\Omega$. Denote by $\mathcal{L}^{1}_{\operatorname{loc}}(\Omega)$ the space of such functions and by $L^{1}_{\operatorname{loc}}(\Omega)$ the space of their equivalence classes modulo equality almost everywhere.
\end{definition}

The following consequence states that a locally integrable function is determined, up to a null set, by its integrals against test functions. It will be used repeatedly to justify uniqueness of weak derivatives.

\begin{proposition}\label{prop unicidad derivada debil prop 14.49}
Let $\Omega\subseteq\mathbb{R}^{n}$ be open. If $f\in L^{1}_{\operatorname{loc}}(\Omega)$ satisfies $\displaystyle\int_{\Omega}f(x)\varphi(x)\,d\lambda_{n}(x)=0$ for every $\varphi\in C_{c}^{\infty}(\Omega)$, then $f=0$ almost everywhere on $\Omega$.
\end{proposition}
\begin{proof}
Suppose, for a contradiction, that $f$ does not vanish almost everywhere.
An exhaustion of $\Omega$ by compact sets and the decomposition
$\displaystyle \{|f|>0\}=\bigcup_{j\geq1}\{|f|\geq j^{-1}\}$ provide a compact set
$K\Subset\Omega$ and a $j\geq1$ for which
$K\cap\{|f|\geq j^{-1}\}$ has positive measure. Divide the unit
circle into finitely many sectors of angular width less than $\displaystyle\frac{\pi}{2}$.
One of them determines a measurable set $A\subseteq K$ of positive
measure and a number $\omega\in\mathbb C$, $|\omega|=1$, such that
\[
 \operatorname{Re}(\omega f)\geq c|f|\geq \frac{c}{j}
 \quad\text{on }A
\]
for some constant $c>0$. In the real case, it suffices to take one of the two
signs.

By regularity of Lebesgue measure, choose a compact set
$C\subseteq A$ of positive measure. Then
$\displaystyle I_C:=\int_C|f|\,d\lambda_n>0$. Absolute continuity of the integral and
outer regularity allow us to choose an open set $U$ with
$C\subseteq U\Subset\Omega$ and
\[
 \int_{U\setminus C}|f|\,d\lambda_n<\frac{c}{2}I_C.
\]
By Proposition~\ref{funciones flan}, there exists
$\eta\in C_c^\infty(U)$ such that $0\leq\eta\leq1$ and $\eta=1$ on $C$.
Applying the hypothesis to $\varphi=\omega\eta$ instead gives
\[
 \operatorname{Re}\int_\Omega f\varphi\,d\lambda_n
 \geq cI_C-\int_{U\setminus C}|f|\,d\lambda_n
 >\frac{c}{2}I_C>0,
\]
a contradiction. Therefore $f=0$ almost everywhere on $\Omega$.
\end{proof}

\subsection{The Bochner integral}

Throughout this subsection, we assume that $(M,\mathcal{A},\mu)$ is complete. If the original measure is not complete, we work with its completion, obtained by adjoining subsets of measurable null sets and extending the measure by assigning them measure zero. We consider functions defined on $(M,\mathcal{A},\mu)$ and taking values in a Banach space $E$. The notion of simple function introduced above remains valid upon taking $V=E$; strong measurability is obtained by pointwise approximation by these functions, whereas the Bochner integral is obtained by approximation in the $L^{1}(M,\mu,E)$ norm.

\begin{definition}\label{def:b5-integral-bochner}\index{strongly measurable function}\index{Bochner integral}
Let $E$ be a Banach space.
\begin{enumerate}[label=(\alph*)]
\item A function $f\colon M\longrightarrow E$ is \textbf{strongly measurable} if there exists a sequence $(s_{k})_{k\in\mathbb{N}}$ of simple functions $s_{k}\colon M\longrightarrow E$ such that $s_{k}(x)\longrightarrow f(x)$ for $\mu$-almost every $x\in M$.
\item A simple function $s\colon M\longrightarrow E$ is \textbf{integrable} if $\displaystyle\int_{M}\|s\|_{E}\,d\mu<\infty$. If $\displaystyle s=\displaystyle\sum_{i=1}^{r}v_{i}\mathbf{1}_{A_{i}}$, where the $A_{i}$ are measurable and pairwise disjoint and terms with $v_{i}=0$ have been omitted, then $\mu(A_{i})<\infty$ for each $i$, and we define $\displaystyle\int_{M}s\,d\mu:=\displaystyle\sum_{i=1}^{r}v_{i}\mu(A_{i})$.
\item A strongly measurable function $f\colon M\longrightarrow E$ is \textbf{Bochner integrable} if there exists a sequence $(s_{k})_{k\in\mathbb{N}}$ of integrable simple functions such that $\displaystyle\int_{M}\|f-s_{k}\|_{E}\,d\mu\longrightarrow0$. In this case, define $\displaystyle\int_{M}f\,d\mu:=\lim_{k\to\infty}\int_{M}s_{k}\,d\mu$. Denote by $\mathcal{L}^{1}(M,\mu,E)$ the space of Bochner integrable functions and by $L^{1}(M,\mu,E)$ the space of their equivalence classes modulo equality almost everywhere. For $A\in\mathcal{A}$, define $\displaystyle\int_{A}f\,d\mu:=\int_{M}f\mathbf{1}_{A}\,d\mu$.
\end{enumerate}
\end{definition}

\begin{lemma}\label{lem:b5-integral-funcion-simple-banach}
The integral of an integrable simple function $s\colon M\longrightarrow E$ is well defined, that is, independent of the representation of $s$, and satisfies $\displaystyle\left\|\int_{M}s\,d\mu\right\|_{E}\leq\int_{M}\|s\|_{E}\,d\mu$.
\end{lemma}
\begin{proof}
Two representations of $s$ admit a common refinement consisting of intersections of their level sets, and the sum defining the integral agrees when computed on this refinement. If $\displaystyle s=\displaystyle\sum_{i=1}^{r}v_{i}\mathbf{1}_{A_{i}}$ with the $A_{i}$ pairwise disjoint, then $\displaystyle\left\|\int_{M}s\,d\mu\right\|_{E}=\left\|\displaystyle\sum_{i=1}^{r}v_{i}\mu(A_{i})\right\|_{E}\leq\displaystyle\sum_{i=1}^{r}\|v_{i}\|_{E}\mu(A_{i})=\int_{M}\|s\|_{E}\,d\mu$.
\end{proof}

The following criterion shows that, once strong measurability has been established, Bochner integrability is characterized entirely by integrability of the norm.

\begin{proposition}\label{prop:b5-criterio-integrabilidad-bochner}
Let $E$ be a Banach space, and let $f\colon M\longrightarrow E$ be strongly measurable. Then $f$ is Bochner integrable if and only if $\|f\|_{E}\in\mathcal{L}^{1}(M,\mu)$. In this case, the Bochner integral is well defined, and the map $\displaystyle\mathcal{I}_{B}\colon \mathcal{L}^{1}(M,\mu,E)\longrightarrow E$, given by $\displaystyle\mathcal{I}_{B}(f):=\int_{M}f\,d\mu$, is linear and satisfies $\displaystyle\left\|\int_{M}f\,d\mu\right\|_{E}\leq\int_{M}\|f\|_{E}\,d\mu$.
\end{proposition}
\begin{proof}
If $(t_k)_{k\in\mathbb N}$ is a sequence of simple functions converging to $f$ almost everywhere, continuity of the norm gives $\|t_k\|_E\to\|f\|_E$ almost everywhere. Since $\mu$ is complete, $\|f\|_E$ is measurable. First suppose that $f$ is Bochner integrable, and let $(s_{k})_{k\in\mathbb{N}}$ be a sequence of integrable simple functions such that $\displaystyle\int_{M}\|f-s_{k}\|_{E}\,d\mu\longrightarrow0$. By Lemma~\ref{lem:b5-integral-funcion-simple-banach}, $\displaystyle\left\|\int_{M}s_{k}\,d\mu-\int_{M}s_{\ell}\,d\mu\right\|_{E}\leq\int_{M}\|s_{k}-s_{\ell}\|_{E}\,d\mu\leq\int_{M}\|s_{k}-f\|_{E}\,d\mu+\int_{M}\|f-s_{\ell}\|_{E}\,d\mu$, so the sequence of integrals is Cauchy. The same argument applied simultaneously to two approximating sequences shows that the limit is independent of the choice of sequence. Moreover, for any $k$ we have $\displaystyle\int_{M}\|f\|_{E}\,d\mu\leq\int_{M}\|f-s_{k}\|_{E}\,d\mu+\int_{M}\|s_{k}\|_{E}\,d\mu<\infty$, and $\displaystyle\left|\int_{M}\|s_{k}\|_{E}\,d\mu-\int_{M}\|f\|_{E}\,d\mu\right|\leq\int_{M}\|s_{k}-f\|_{E}\,d\mu\longrightarrow0$. Consequently, $\displaystyle\left\|\int_{M}f\,d\mu\right\|_{E}=\lim_{k\to\infty}\left\|\int_{M}s_{k}\,d\mu\right\|_{E}\leq\lim_{k\to\infty}\int_{M}\|s_{k}\|_{E}\,d\mu=\int_{M}\|f\|_{E}\,d\mu$. Linearity follows by approximating two Bochner integrable functions by simple functions and taking linear combinations.

Conversely, suppose that $\|f\|_{E}\in\mathcal{L}^{1}(M,\mu)$. By strong measurability, there exists a sequence of simple functions $(t_{k})_{k\in\mathbb{N}}$ such that $t_{k}(x)\longrightarrow f(x)$ for $\mu$-almost every $x\in M$. Define $B_{k}:=\{x\in M\mid \|t_{k}(x)\|_{E}\leq2\|f(x)\|_{E}\}$ and $s_{k}:=t_{k}\mathbf{1}_{B_{k}}$. Each $B_{k}$ is measurable, $s_{k}$ is simple, and $\|s_{k}\|_{E}\leq2\|f\|_{E}$, so $s_{k}$ is integrable. If $f(x)\neq0$, then $\|t_{k}(x)\|_{E}\longrightarrow\|f(x)\|_{E}$, and hence there exists $k_x\in\mathbb N$ such that $x\in B_k$ for every $k\geq k_x$; if $f(x)=0$, then $s_k(x)=0$ for every $k\in\mathbb N$. Thus $s_{k}(x)\longrightarrow f(x)$ for $\mu$-almost every $x\in M$. Finally, $\|s_{k}-f\|_{E}\leq3\|f\|_{E}$ and Theorem~\ref{convergencia dominada} imply that $\displaystyle\int_{M}\|s_{k}-f\|_{E}\,d\mu\longrightarrow0$. Consequently, $f$ is Bochner integrable.
\end{proof}

Stability of strong measurability under almost everywhere pointwise limits will be used in the dominated convergence theorem.

\begin{lemma}\label{lem:b5-limite-fuertemente-medible}
Let $(f_{k})_{k\in\mathbb{N}}$ be a sequence of strongly measurable functions $f_{k}\colon M\longrightarrow E$. If $f_{k}(x)\longrightarrow f(x)$ for $\mu$-almost every $x\in M$, then $f$ is strongly measurable.
\end{lemma}
\begin{proof}
For each $k$, the essential range of $f_{k}$ is contained in a separable subspace of $E$, since $f_{k}$ is an almost everywhere pointwise limit of simple functions. Thus, outside a null set, the ranges of all the $f_{k}$ and of $f$ lie in a single closed separable subspace $E_{0}\subseteq E$. Moreover, each $f_{k}$ is measurable and, since $E$ is metrizable, its almost everywhere pointwise limit $f$ is also measurable after modification on a null set. Let $(e_{j})_{j\in\mathbb{N}}$ be a dense subset of $E_{0}$. For each $r\in\mathbb{N}$, choose at every point an element of $\{e_{1},\dots,e_{r}\}$ at minimum distance from $f(x)$, resolving ties by choosing the smallest index. The resulting functions are simple and converge pointwise to $f$, since $\displaystyle\min_{1\leq j\leq r}\|f(x)-e_{j}\|_{E}\longrightarrow0$. Consequently, $f$ is strongly measurable.
\end{proof}

The vector-valued dominated convergence theorem is now a direct consequence of the preceding criterion.

\begin{theorem}[Bochner dominated convergence theorem]\label{teo:b5-convergencia-dominada-bochner}\index{convergence theorem!Bochner dominated}
Let $(f_{k})_{k\in\mathbb{N}}$ be a sequence of strongly measurable functions $f_{k}\colon M\longrightarrow E$ such that $f_{k}(x)\longrightarrow f(x)$ for $\mu$-almost every $x\in M$. Suppose that there exists $g\in\mathcal{L}^{1}(M,\mu)$ such that $\|f_{k}\|_{E}\leq g$ $\mu$-almost everywhere for every $k\in\mathbb{N}$. Then $f$ is strongly measurable, $f_{k},f\in\mathcal{L}^{1}(M,\mu,E)$, and $\displaystyle\int_{M}\|f_{k}-f\|_{E}\,d\mu\longrightarrow0$. In particular, $\displaystyle\int_{M}f_{k}\,d\mu\longrightarrow\int_{M}f\,d\mu$ in $E$.
\end{theorem}
\begin{proof}
By Lemma~\ref{lem:b5-limite-fuertemente-medible}, $f$ is strongly measurable. Moreover, $\|f\|_{E}\leq g$ almost everywhere and $\|f_{k}-f\|_{E}\leq2g$. Theorem~\ref{convergencia dominada} applied to the scalar functions $\|f_{k}-f\|_{E}$ implies convergence in $\mathcal{L}^{1}$, while Proposition~\ref{prop:b5-criterio-integrabilidad-bochner} ensures Bochner integrability. Finally, $\displaystyle\left\|\int_{M}f_{k}\,d\mu-\int_{M}f\,d\mu\right\|_{E}\leq\int_{M}\|f_{k}-f\|_{E}\,d\mu\longrightarrow0$.
\end{proof}

\begin{theorem}[Fubini's theorem for Bochner integrals]
\label{teo:b5-fubini-bochner}
\index{Fubini's theorem!for Bochner integrals}
Let $(X,\mathcal A,\mu)$ and $(Y,\mathcal B,\nu)$ be
complete $\sigma$-finite measure spaces, let $E$ be a Banach space, and let
$F\colon X\times Y\longrightarrow E$ be strongly measurable with respect to the completion $\overline{\mu\otimes\nu}$ of the product measure. If
\[
\int_{X\times Y}\|F(x,y)\|_E\,d\overline{\mu\otimes\nu}(x,y)<\infty,
\]
then $F(x,\cdot)$ is Bochner integrable for $\mu$-almost every $x$, and
$F(\cdot,y)$ is Bochner integrable for $\nu$-almost every $y$. The functions
obtained by integrating in one variable are strongly measurable and integrable in
the other variable, and
\begin{align*}
\int_{X\times Y}F\,d\overline{\mu\otimes\nu}
&=
\int_X\left(\int_YF(x,y)\,d\nu(y)\right)d\mu(x)\\
&=
\int_Y\left(\int_XF(x,y)\,d\mu(x)\right)d\nu(y).
\end{align*}
\end{theorem}

\begin{proof}
The criterion in Proposition~\ref{prop:b5-criterio-integrabilidad-bochner}
shows that $F$ is Bochner integrable on the completed product.
A sequence of simple functions converges to $F$ outside a product
null set. The completed version of Tonelli,
Theorem~\ref{teo:tonelli}, allows us to discard a $\mu$-null set
of values of $x$ so that all these simple sections are measurable
and converge to $F(x,\cdot)$ $\nu$-almost everywhere.
Furthermore, Tonelli applied to $\|F\|_E$ gives
\[
 \int_X\left(\int_Y\|F(x,y)\|_E\,d\nu(y)\right)d\mu(x)<\infty.
\]
Thus almost every section is strongly measurable and has integrable
norm; the same criterion proves its Bochner integrability.

Choose integrable simple functions $S_k$ such that
\[
 \int_{X\times Y}\|F-S_k\|_E\,d\overline{\mu\otimes\nu}<2^{-k}.
\]
To describe the partial integrals of a simple function, write
$S=\sum_{i=1}^l v_i\mathbf1_{A_i}$, with $v_i\neq0$ and disjoint
measurable sets in the completed product.
Integrability of $S$ implies
$\overline{\mu\otimes\nu}(A_i)<\infty$.
Tonelli shows that $a_i(x):=\nu((A_i)_x)$ are measurable,
integrable and finite $\mu$-almost everywhere.
Defining them to be zero on the exceptional sets, we obtain
\[
 G_S(x):=\int_Y S(x,y)\,d\nu(y)=\sum_{i=1}^l v_i a_i(x)
\]
almost everywhere. Approximating each scalar $a_i$ by simple
functions shows that $G_S$ is strongly measurable.
The inequality for the integral of a simple function and Tonelli give
\[
 \int_X\|G_S(x)\|_E\,d\mu(x)
 \leq\int_{X\times Y}\|S\|_E\,d\overline{\mu\otimes\nu}<\infty.
\]
Hence $G_S$ is Bochner integrable. Linearity and scalar Tonelli
applied to each $\mathbf1_{A_i}$ prove
\[
 \int_X G_S\,d\mu
 =\sum_{i=1}^l v_i\overline{\mu\otimes\nu}(A_i)
 =\int_{X\times Y}S\,d\overline{\mu\otimes\nu}.
\]

Put $G_k:=G_{S_k}$ and define
$G(x):=\int_YF(x,y)\,d\nu(y)$ where the section is integrable,
and $G(x)=0$ on the measurable exceptional set. By Tonelli,
\[
 \int_X\sum_{k\geq1}
      \left(\int_Y\|F(x,y)-S_k(x,y)\|_E\,d\nu(y)\right)d\mu(x)
 <\sum_{k\geq1}2^{-k}<\infty.
\]
The inner series is finite almost everywhere, so its terms tend
to zero. The bound for the Bochner integral then gives
$G_k(x)\to G(x)$ almost everywhere.
Lemma~\ref{lem:b5-limite-fuertemente-medible} shows that $G$ is
strongly measurable. Since
$\|G(x)\|_E\leq\int_Y\|F(x,y)\|_E\,d\nu(y)$ almost everywhere,
$G$ is Bochner integrable. Moreover,
\[
 \int_X\|G-G_k\|_E\,d\mu
 \leq\int_{X\times Y}\|F-S_k\|_E\,d\overline{\mu\otimes\nu}
 <2^{-k}.
\]
We may therefore pass to the limit in the identity for $S_k$
to obtain the first equality. Interchanging $X$ and $Y$
proves all the claims for the other sections and the second equality.
\end{proof}

Like the Riemann integral, the Bochner integral is compatible with continuous and closed linear operators.

\begin{proposition}\label{prop:b5-integral-bochner-operadores-lineales}
Let $E$ and $F$ be Banach spaces, and let $f\colon M\longrightarrow E$ be Bochner integrable.
\begin{enumerate}[label=(\alph*)]
\item If $T\in\mathcal{L}(E,F)$, then $T\circ f$ is Bochner integrable and $\displaystyle\int_{M}T(f)\,d\mu=T\left(\int_{M}f\,d\mu\right)$.
\item If $T\colon \operatorname{Dom}(T)\subseteq E\longrightarrow F$ is a closed linear operator, $f(x)\in\operatorname{Dom}(T)$ for $\mu$-almost every $x\in M$, and $T\circ f$ is Bochner integrable, then $\displaystyle\int_{M}f\,d\mu\in\operatorname{Dom}(T)$ and the same equality holds.
\end{enumerate}
\end{proposition}
\begin{proof}
If $T\in\mathcal{L}(E,F)$ and $(s_{k})_{k\in\mathbb{N}}$ approximates $f$ in $\mathcal{L}^{1}$, then $T\circ s_{k}$ is simple, $\displaystyle\int_{M}\|T(f)-T(s_{k})\|_{F}\,d\mu\leq\|T\|\int_{M}\|f-s_{k}\|_{E}\,d\mu\longrightarrow0$, and the equality follows first for simple functions and then by passage to the limit.

For the second assertion, modifying $f$ and $T\circ f$ on a null set, we may assume that $f(x)\in\operatorname{Dom}(T)$ for every $x\in M$. Consider $G\colon M\longrightarrow E\times F$, given by $G(x):=(f(x),T(f(x)))$, and equip $E\times F$ with the norm $\|(u,v)\|:=\|u\|_{E}+\|v\|_{F}$. The function $G$ is strongly measurable and $\displaystyle\int_{M}\|G\|\,d\mu<\infty$, so it is Bochner integrable. The coordinate projections and the first part show that $\displaystyle\int_{M}G\,d\mu=\left(\int_{M}f\,d\mu,\int_{M}T(f)\,d\mu\right)$. Since the graph $\mathcal{G}(T)$ is a closed subspace of $E\times F$, the quotient map $Q\colon E\times F\longrightarrow(E\times F)/\mathcal{G}(T)$ is linear and continuous. Because $Q\circ G=0$ almost everywhere, the first part implies $\displaystyle Q\left(\int_{M}G\,d\mu\right)=\int_{M}Q(G)\,d\mu=0$. Thus $\displaystyle\int_{M}G\,d\mu\in\mathcal{G}(T)$, which is equivalent to the conclusion.
\end{proof}

\begin{remark}\label{obs:b5-integral-bochner-pettis}
If $f\colon M\longrightarrow E$ is Bochner integrable and $\varphi\in E'$, Proposition~\ref{prop:b5-integral-bochner-operadores-lineales} applied to $\varphi\colon E\longrightarrow\mathbb{K}$ shows that $\varphi\circ f$ is integrable and $\displaystyle\varphi\left(\int_{A}f\,d\mu\right)=\int_{A}\varphi(f)\,d\mu$ for every $A\in\mathcal{A}$. In particular, every Bochner integrable function is Pettis integrable and hence Dunford integrable; all three integrals agree under the canonical inclusion $E\hookrightarrow E''$.
\end{remark}

\begin{theorem}[Radon--Nikodym]\label{teo:radon-nikodym}\index{Radon--Nikodym theorem}
Let $\mu$ and $\nu$ be $\sigma$-finite measures on $(X,\mathcal A)$, with $\nu$ positive and $\nu\ll\mu$. There exists a measurable function $h\colon X\longrightarrow[0,\infty)$, unique $\mu$-almost everywhere, such that
\[
\nu(A)=\int_Ah\,d\mu,
\qquad A\in\mathcal A.
\]
If $\nu$ is a signed or complex measure, $\nu\ll\mu$, and its total
variation $|\nu|$ is $\sigma$-finite, there similarly exists a measurable
density $h\colon X\to\mathbb K$. More precisely, one may choose a
increasing sequence of measurable sets $(X_k)_{k\in\mathbb N}$, with
$\displaystyle X=\bigcup_{k=1}^{\infty}X_k$ and
$\mu(X_k)+|\nu|(X_k)<\infty$, and then
$h\mathbf 1_{X_k}\in L^1(X,\mu)$ and
\[
 \nu(A\cap X_k)=\int_{A\cap X_k}h\,d\mu
 \qquad(A\in\mathcal A,\ k\in\mathbb N).
\]
In particular, if $|\nu|(X)<\infty$, we have $h\in L^1(X,\mu)$. This
formulation replaces the notation $L^1_{\mathrm{loc}}(X,\mu)$, which is not
canonical on an abstract measure space without a topological notion of
locality.
\end{theorem}
\begin{proof}
First suppose that $\mu(X)+\nu(X)<\infty$ and put
$\lambda:=\mu+\nu$. Work in the real Hilbert space
$L^2(X,\lambda;\mathbb R)$. The functional
\[
 L(f):=\int_Xf\,d\nu
\]
is well defined and continuous, since
$|L(f)|\leq\nu(X)^{\displaystyle\frac{1}{2}}\|f\|_{L^2(X,\lambda)}$.
Theorem~\ref{teo:representacion-riesz-hilbert} gives
$g\in L^2(X,\lambda;\mathbb R)$ such that
\[
 \nu(A)=\int_Ag\,d\lambda,
 \qquad
 \mu(A)=\int_A(1-g)\,d\lambda,
 \qquad A\in\mathcal A.
\]
Positivity of the two measures implies $0\leq g\leq1$
$\lambda$-almost everywhere. Indeed, applying these identities
to $\{g\leq-\displaystyle\frac{1}{m}\}$ and $\{g\geq1+\displaystyle\frac{1}{m}\}$, for $m\in\mathbb N$,
shows that all these sets are null. Choose a measurable representative
with values in $[0,1]$. On $Z:=\{g=1\}$ the measure $\mu$ vanishes;
since $\nu\ll\mu$, so does $\nu$, and $\lambda(Z)=0$.
Define $h=\displaystyle\frac{g}{1-g}$ on $X\setminus Z$ and $h=0$ on $Z$.

The identity $\mu(A)=\int_A(1-g)\,d\lambda$ implies
\[
 \int_Xu\,d\mu=\int_Xu(1-g)\,d\lambda
 \qquad\text{for every measurable function }u\geq0.
\]
This is first verified for indicators, then for simple functions,
and finally by monotone convergence. Applying it to $u=h\mathbf1_A$
gives $\int_Ah\,d\mu=\int_Ag\,d\lambda=\nu(A)$.

Prove uniqueness when $\nu(X)<\infty$. If $h_1,h_2\geq0$
represent $\nu$, both are integrable. For
$E_m:=\{h_1\geq h_2+\displaystyle\frac{1}{m}\}$ we have
\[
 0=\int_{E_m}(h_1-h_2)\,d\mu
 \geq\frac{1}{m}\mu(E_m).
\]
All integrals in this equality are finite. Thus,
$\mu\{h_1>h_2\}=0$; interchanging the two functions proves uniqueness.

In the $\sigma$-finite case, choose increasing measurable sets
$X_k$ covering $X$ with $\mu(X_k)+\nu(X_k)<\infty$.
Put $X_0=\varnothing$ and $D_k=X_k\setminus X_{k-1}$.
Apply the finite case on each $D_k$ and define $h$ using
the resulting densities on this measurable partition. Countable
additivity and monotone convergence give
\[
 \int_Ah\,d\mu
 =\sum_{k=1}^{\infty}\int_{A\cap D_k}h\,d\mu
 =\sum_{k=1}^{\infty}\nu(A\cap D_k)=\nu(A).
\]
If two densities represent $\nu$, their restrictions to each $X_k$
are integrable and agree by finite-case uniqueness. The countable union
of the exceptional sets proves uniqueness on $X$.

For a signed or complex measure, now choose an increasing exhaustion
with $\mu(X_k)+|\nu|(X_k)<\infty$. On each $X_k$ apply the positive
case to the Jordan parts of the real component and, in the complex
case, of the imaginary component. This gives a density
$h_k\in L^1(X_k,\mu)$ for $\nu\restriction_{X_k}$.
Uniqueness for integrable real densities follows by applying
the $E_m$ argument to their difference; for complex densities apply
it to both components. Hence $h_{k+1}=h_k$ almost everywhere
on $X_k$. Defining $h$ by $h_k$ on $D_k$ gives a measurable function
that agrees almost everywhere with each $h_k$ on $X_k$ and satisfies
the stated identities. Moreover,
\[
 \int_{X_k}|h|\,d\mu
 \leq |\operatorname{Re}\nu|(X_k)
       +|\operatorname{Im}\nu|(X_k)
 \leq2|\nu|(X_k),
\]
with the imaginary component omitted in the real case. If $|\nu|(X)<\infty$,
this bound and monotone convergence show that $h\in L^1(X,\mu)$.
Global uniqueness again follows on the sets $X_k$.
\end{proof}

\section{The \texorpdfstring{$L^{p}(X)$}{Lp} spaces}\label{ap:espacios-Lp}

The $L^p(X)$ spaces will be the starting point for Sobolev spaces and several weak formulations of partial differential equations. We recall their basic properties and some inequalities used later.

\begin{definition}\label{def:b6-espacios-lp-norma-espacio-medida-fijo}\index{Lebesgue space!integrable measurable functions}\index{norm!on Lebesgue spaces}
Let $(X,\mathcal{A},\mu)$ be a measure space, and fix $p\in (0,\infty)$. Define \[\mathcal{L}^{p}(X):=\left\{f\colon X\longrightarrow\mathbb{K}\;\middle|\;\text{$f$ is measurable and }\int_{X}|f|^{p}d\mu<\infty\right\}.\] Define $\|\cdot\|_{\mathcal L^p(X)}\colon \mathcal{L}^{p}(X)\longrightarrow\mathbb{R}$ by
\[
 \|f\|_{\mathcal L^p(X)}
 =\left(\displaystyle\int_{X}|f|^{p}\,d\mu\right)^{\frac{1}{p}}.
 \]

If $p=\infty$, define \[\mathcal{L}^{\infty}(X):=\left\{f\colon X\longrightarrow\mathbb{K}\mid \text{$f$ is measurable and bounded almost everywhere on $X$}\right\},\] and
\[
 \|f\|_{\mathcal L^\infty(X)}
 :=\inf\{a>0\mid\mu(\{x\in X\mid |f(x)|>a\})=0\}.
 \]
\end{definition}
\begin{definition}\label{def:b6-espacios-lp-norma-fijo-relacion-sigue-note}\index{Lebesgue space as a quotient space}
Fix $p\in[1,\infty]$. Define a relation $\sim$ on $\mathcal{L}^{p}(X)$ as follows:
\[f\sim g\iff \|f-g\|_{\mathcal L^p(X)}=0.\] Note that $\sim$ is an equivalence relation; write $L^{p}(X):=\mathcal{L}^{p}(X)/\sim$. If $[f]\in L^{p}(X)$, define $\|[f]\|_{L^p(X)}=\|f\|_{\mathcal L^p(X)}$, and in practice identify $[f]$ with $f$.
\end{definition}

\begin{theorem}[Lebesgue spaces]\label{teo:b6-espacios-lp-espacios-de-lebesgue}\index{Lebesgue space} Let $p\in[1,\infty]$. Then:
\begin{enumerate}[label=(\alph*)]
\item If $p<\infty$,
$\|\cdot\|_{L^p(X)}$ is a norm on $L^{p}(X)$ and $(L^{p}(X),\|\cdot\|_{L^p(X)})$ is a Banach space. The triangle inequality in $L^{p}(X)$ is known as Minkowski's inequality.
\item If $p=\infty$, $\|\cdot\|_{L^\infty(X)}$ is a norm on $L^{\infty}(X)$ and $(L^{\infty}(X),\|\cdot\|_{L^\infty(X)})$ is a Banach space.
\end{enumerate}

\end{theorem}
\begin{proof}
First consider \(1\leq p<\infty\). Homogeneity and separation of points modulo equality almost everywhere follow directly from the definition. To prove the triangle inequality when \(1<p<\infty\), first recall Hölder's inequality. If \(u\in L^p(X)\), \(v\in L^{p'}(X)\), and both norms are nonzero, Young's inequality applied to
\[
\frac{|u|}{\|u\|_{L^{p}(X)}}
\quad\text{and}\quad
\frac{|v|}{\|v\|_{L^{p'}(X)}}
\]
and integrated over \(X\) gives
\[
\int_X|uv|\,d\mu\leq\|u\|_{L^{p}(X)}\|v\|_{L^{p'}(X)}.
\]
If either norm vanishes, the same inequality is immediate.
Before applying it to $f+g$, verify that this function belongs to
$L^p(X)$ without using Minkowski yet. Convexity of $t\mapsto t^p$
on $[0,\infty)$ gives
\[
 |f+g|^p\leq(|f|+|g|)^p
 \leq2^{p-1}(|f|^p+|g|^p).
\]
The right-hand side is integrable; hence $f+g\in L^p(X)$ and
$|f+g|^{p-1}\in L^{p'}(X)$ because $(p-1)p'=p$.
Applying Hölder with \(v=|f+g|^{p-1}\) gives
\[
\begin{aligned}
\|f+g\|_{L^{p}(X)}^p
&\leq
\int_X|f|\,|f+g|^{p-1}\,d\mu
+
\int_X|g|\,|f+g|^{p-1}\,d\mu\\
&\leq
\bigl(\|f\|_{L^{p}(X)}+\|g\|_{L^{p}(X)}\bigr)
\bigl\||f+g|^{p-1}\bigr\|_{L^{p'}(X)}\\
&=
\bigl(\|f\|_{L^{p}(X)}+\|g\|_{L^{p}(X)}\bigr)\|f+g\|_{L^{p}(X)}^{p-1}.
\end{aligned}
\]
If \(\|f+g\|_{L^{p}(X)}=0\), there is nothing to prove; otherwise, divide by \(\|f+g\|_{L^{p}(X)}^{p-1}\). For \(p=1\), the inequality follows by integrating \(|f+g|\leq|f|+|g|\). This proves that \(\|\cdot\|_{L^p(X)}\) is a norm.

We now prove completeness. Let \((f_n)\) be a Cauchy sequence in \(L^p(X)\). We may choose a subsequence \((f_{n_k})\) such that
\[
\|f_{n_{k+1}}-f_{n_k}\|_{L^{p}(X)}\leq2^{-k},
\qquad k\in\mathbb N.
\]
Set \(g_k:=|f_{n_{k+1}}-f_{n_k}|\). By Minkowski's inequality, already proved,
\[
\left\|\sum_{k=1}^Ng_k\right\|_{L^{p}(X)}
\leq\sum_{k=1}^N\|g_k\|_{L^{p}(X)}
\leq1
\qquad(N\in\mathbb N).
\]
The monotone convergence theorem applied to \(\displaystyle (\displaystyle\sum_{k=1}^Ng_k)^p\) shows that
\[
G:=\sum_{k=1}^{\infty}g_k
\quad\text{belongs to }L^p(X)
\quad\text{and}\quad
\|G\|_{L^{p}(X)}\leq1.
\]
Outside a null set, the series
\[
\sum_{k=1}^{\infty}\bigl(f_{n_{k+1}}-f_{n_k}\bigr)
\]
converges absolutely. Define there
\[
f:=f_{n_1}+
\sum_{k=1}^{\infty}\bigl(f_{n_{k+1}}-f_{n_k}\bigr)
\]
and define $f=0$ on the measurable exceptional set.
The function $f$ is measurable, being a limit of measurable functions
on its complement and constant on that set. This choice preserves
measurability even when the measure is not complete. Moreover, for each \(j\),
\[
|f-f_{n_j}|
\leq\sum_{k=j}^{\infty}g_k
\]
almost everywhere. Applying monotone convergence again to the tails and using Minkowski for their finite sums gives
\[
\|f-f_{n_j}\|_{L^{p}(X)}
\leq\sum_{k=j}^{\infty}2^{-k}\longrightarrow0.
\]
In particular, \(f\in L^p(X)\). Since the original sequence is Cauchy and a subsequence converges to \(f\), the entire sequence converges to \(f\).

Now consider \(p=\infty\). The norm properties follow from the pointwise inequalities
\[
|\lambda f|=|\lambda|\,|f|,
\qquad
|f+g|\leq|f|+|g|
\]
outside null sets. If \((f_n)\) is Cauchy in \(L^\infty(X)\), choose a subsequence such that
\[
\|f_{n_{k+1}}-f_{n_k}\|_\infty\leq2^{-k}.
\]
After modifying representatives on the union of a countable family of null sets, we may assume that
\[
|f_{n_{k+1}}(x)-f_{n_k}(x)|\leq2^{-k}
\]
for every \(k\) and every \(x\) outside a single null set \(N\). By the Weierstrass test, \(f_{n_k}\) converges uniformly on \(X\setminus N\) to a measurable function \(f\), which we define to be zero on \(N\), and
\[
\|f-f_{n_j}\|_\infty
\leq\sum_{k=j}^{\infty}2^{-k}\longrightarrow0.
\]
Again, the Cauchy property implies that the entire sequence converges to \(f\) in \(L^\infty(X)\). Thus \(L^\infty(X)\) is complete as well.
\end{proof}

\begin{proposition}[Reverse Minkowski inequality for exponents less than one]
\label{prop:minkowski-inversa-exponentes-menores-uno}
Let $(X,\mathcal A,\mu)$ be a measure space, let $0<r<1$, and let
$f,g\colon X\longrightarrow[0,\infty]$ be measurable functions such that
$\displaystyle\int_Xf^r\,d\mu<\infty$ and
$\displaystyle\int_Xg^r\,d\mu<\infty$. Then
\[
\left(\int_X(f+g)^r\,d\mu\right)^{\frac{1}{r}}
\geq
\left(\int_Xf^r\,d\mu\right)^{\frac{1}{r}}
+
\left(\int_Xg^r\,d\mu\right)^{\frac{1}{r}}.
\]
\end{proposition}

\begin{proof}
Set
\[
A:=\left(\int_Xf^r\,d\mu\right)^{\frac{1}{r}},
\qquad
B:=\left(\int_Xg^r\,d\mu\right)^{\frac{1}{r}}.
\]
If $A=0$ or $B=0$, one of the functions vanishes almost everywhere and the
inequality is immediate. Thus suppose that $A,B>0$ and define
$\theta:=\displaystyle\frac{A}{A+B}$. Concavity of
$t\mapsto t^r$ on $[0,\infty)$ gives
\[
\left(\theta\frac{f}{A}+(1-\theta)\frac{g}{B}\right)^r
\geq
\theta\left(\frac{f}{A}\right)^r
+(1-\theta)\left(\frac{g}{B}\right)^r.
\]
Upon integration, the right-hand side has integral one. Since
\[
f+g
=
(A+B)\left(\theta\frac{f}{A}+(1-\theta)\frac{g}{B}\right),
\]
we obtain
\[
\int_X(f+g)^r\,d\mu\geq(A+B)^r.
\]
Taking the $r$th root completes the proof.
\end{proof}

\begin{proposition}[Density of simple functions]
\label{prop:densidad-funciones-simples-Lp}
Let $(X,\mathcal A,\mu)$ be a $\sigma$-finite measure space, and let
$1\leq p<\infty$. Simple functions whose support has finite measure
are dense in $L^p(X)$.
\end{proposition}

\begin{proof}
Choose increasing sets $X_k\in\mathcal A$ of finite measure with
$\displaystyle\bigcup_{k\in\mathbb N}X_k=X$. For $f\in L^p(X)$, the
dominated convergence theorem~\ref{convergencia dominada}
gives
\[
f\mathbf1_{X_k\cap\{|f|\leq k\}}\longrightarrow f
\quad\text{in }L^p(X).
\]
Each function on the left-hand side is bounded and has support of finite
measure. Approximate its positive and negative parts from below by
step functions with step size $2^{-m}$; the pointwise error is at most
$2^{-m}$ on $X_k$, so the error in $L^p(X)$ tends to zero. For
complex-valued functions, apply the argument to the real and imaginary parts.
\end{proof}

\begin{lemma}[Hölder's inequality for two factors]
\label{lem:holder-dos-factores-previo-dualidad}
Let $1\leq p\leq\infty$, and let $p'$ be its conjugate exponent. If
$f\in L^p(X)$ and $g\in L^{p'}(X)$, then $fg\in L^1(X)$ and
\[
 \int_X|fg|\,d\mu\leq\|f\|_{L^p(X)}\|g\|_{L^{p'}(X)}.
\]
\end{lemma}

\begin{proof}
The cases $p=1$ and $p=\infty$ follow directly from the definition of the
essential supremum. If $1<p<\infty$ and both norms are nonzero, apply
Young's inequality
$ab\leq \displaystyle\frac{a^p}{p}+\displaystyle\frac{b^{p'}}{p'}$ pointwise to
$a=\displaystyle\frac{|f|}{\|f\|_{L^{p}(X)}}$ and $b=\displaystyle\frac{|g|}{\|g\|_{L^{p'}(X)}}$ and integrate. If either norm
vanishes, the conclusion is immediate.
\end{proof}

\begin{theorem}[Duality, reflexivity, and separability of $L^p(X)$]
\label{teo:dualidad-reflexividad-Lp}
Let $(X,\mathcal A,\mu)$ be a $\sigma$-finite measure space, and let $1\leq p<\infty$. If $p'$ is the conjugate exponent, the map
\[
\overline{L^{p'}(X)}\longrightarrow (L^p(X))',
\qquad
\overline g\longmapsto\left[f\mapsto\int_Xf\overline g\,d\mu\right]
\]
is a complex-linear isometric isomorphism (in the real case, the bars are
omitted). Equivalently, the direct map
$\displaystyle g\mapsto[f\mapsto\int_Xf\overline g\,d\mu]$ is antilinear. Consequently,
$L^p(X)$ is reflexive if $1<p<\infty$. If, in addition, the $\sigma$-algebra
modulo null sets is generated by a countable family---in particular,
for an open subset of $\mathbb R^n$ with Lebesgue measure---then $L^p(X)$
is separable.
\end{theorem}
\begin{proof}
Lemma~\ref{lem:holder-dos-factores-previo-dualidad} shows that
$\overline g$ defines a functional of norm at most $\|g\|_{L^{p'}(X)}$.
The map is complex-linear because
$\alpha\overline g=\overline{\overline\alpha g}$ in the conjugate space.
Conversely, let $\Lambda\in(L^p(X))'$ and choose increasing measurable
sets $X_k$ with $\mu(X_k)<\infty$ and $\bigcup_kX_k=X$.
Identify functions on $X_k$ with their extensions by zero.
For measurable $A\subseteq X_k$, put $\nu_k(A):=\Lambda(\mathbf1_A)$.
If $(A_j)$ is a disjoint family, the partial sums of its indicators
converge in $L^p(X_k)$ to the indicator of the union, since $p<\infty$
and $\mu(X_k)<\infty$. Continuity of $\Lambda$ proves countable
additivity of $\nu_k$. For a finite partition $(A_j)_{j=1}^N$ of
a subset of $X_k$, choose $|\theta_j|=1$ with
$\theta_j\nu_k(A_j)=|\nu_k(A_j)|$. Then
\[
 \sum_{j=1}^N|\nu_k(A_j)|
 =\left|\Lambda\left(\sum_{j=1}^N\theta_j\mathbf1_{A_j}\right)\right|
 \leq\|\Lambda\|\,\mu(X_k)^{\frac{1}{p}}.
\]
Thus $\nu_k$ has finite variation and $\nu_k\ll\mu$.
Theorem~\ref{teo:radon-nikodym} gives
$g_k\in L^1(X_k,\mu)$ with $d\nu_k=\overline g_k\,d\mu$.
Uniqueness implies $g_{k+1}=g_k$ almost everywhere on $X_k$.
Defining $g$ by $g_k$ on $X_k\setminus X_{k-1}$, with
$X_0=\varnothing$, gives a measurable function agreeing almost
everywhere with $g_k$ on each $X_k$.

The identity
\[
 \Lambda(f)=\int_Xf\overline g\,d\mu
\]
first holds for simple functions supported in $X_k$.
Every bounded measurable function supported in $X_k$ has uniform
simple approximations. These converge in $L^p(X)$ and,
since $g\mathbf1_{X_k}\in L^1(X)$, also in the preceding integral.
The identity is therefore established for these bounded functions.

For $1<p<\infty$, define $\operatorname{sgn}(z)=\displaystyle\frac{z}{|z|}$ if
$z\neq0$ and $\operatorname{sgn}(0)=0$. Apply the identity to
\[
 f_{N,k}:=\operatorname{sgn}(g)|g|^{p'-1}
             \mathbf1_{X_k\cap\{|g|\leq N\}}.
\]
This function is bounded and supported in $X_k$. If
$A_{N,k}:=\int_{X_k\cap\{|g|\leq N\}}|g|^{p'}\,d\mu$, we obtain
\[
 A_{N,k}\leq\|\Lambda\|A_{N,k}^{\frac{1}{p}},
 \qquad A_{N,k}^{\frac{1}{p'}}\leq\|\Lambda\|,
\]
with the second inequality immediate when $A_{N,k}=0$.
Letting $N\to\infty$ and then $k\to\infty$ by monotone convergence
proves $g\in L^{p'}(X)$ and $\|g\|_{L^{p'}(X)}\leq\|\Lambda\|$.
If $p=1$, for $\varepsilon>0$ put
$B_{k,\varepsilon}:=X_k\cap\{|g|>\|\Lambda\|+\varepsilon\}$ and
test with $f=\operatorname{sgn}(g)\mathbf1_{B_{k,\varepsilon}}$.
We obtain
\[
 (\|\Lambda\|+\varepsilon)\mu(B_{k,\varepsilon})
 \leq\int_{B_{k,\varepsilon}}|g|\,d\mu
 =|\Lambda(f)|\leq\|\Lambda\|\mu(B_{k,\varepsilon}).
\]
Hence $\mu(B_{k,\varepsilon})=0$. Taking $k\in\mathbb N$ and
$\varepsilon=\displaystyle\frac{1}{m}$, $m\in\mathbb N$, gives
$\|g\|_{L^\infty(X)}\leq\|\Lambda\|$.
Hölder and density of simple functions of finite measure support
now extend the representation to every $f\in L^p(X)$.
The initial Hölder bound and the reverse bound just proved give
$\|\Lambda\|=\|g\|_{L^{p'}(X)}$. Radon--Nikodym uniqueness
on each $X_k$ also proves uniqueness of $g$.

To prove reflexivity when $1<p<\infty$, use the bilinear
formulation of the same result:
\[
 S\colon L^{p'}(X)\longrightarrow(L^p(X))',
 \qquad (Sh)(f):=\int_Xfh\,d\mu.
\]
This is a linear isometric isomorphism, obtained by replacing $g$
with $\overline h$ in the preceding representation. If
$\Psi\in(L^p(X))''$, then $\Psi\circ S$ is a continuous linear
functional on $L^{p'}(X)$. The representation already proved, with
exponent $p'$, gives $f\in L^p(X)$ such that
\[
 \Psi(Sh)=\int_Xhf\,d\mu=(Sh)(f)=(Jf)(Sh),
 \qquad h\in L^{p'}(X),
\]
where $J$ is the canonical inclusion into the bidual. Surjectivity
of $S$ implies $\Psi=Jf$; thus $J$ is surjective.

Finally, let $\mathcal G$ be a countable family generating
$\mathcal A$ modulo null sets. Add the sets $X_k$ to $\mathcal G$
and denote the countable algebra they generate by $\mathcal C$.
For each $k$, every measurable set $B\subseteq X_k$ can be
approximated in measure by sets $C\cap X_k$, $C\in\mathcal C$.
Indeed, the family of sets admitting this approximation contains
the trace algebra and is closed under complements in $X_k$
and finite unions. It is also closed under increasing unions:
if $B_j\uparrow B$, then $\mu(B\setminus B_j)\to0$, since
$\mu(X_k)<\infty$, and one then approximates a fixed $B_j$.
These properties make it a $\sigma$-algebra. It also contains
all null sets, approximated by the empty set, so the hypothesis
on $\mathcal G$ includes all measurable sets of $X_k$ in this family.

The identity
\[
 \|\mathbf1_B-\mathbf1_{C\cap X_k}\|_{L^p(X)}
 =\mu\bigl(B\mathbin{\triangle}(C\cap X_k)\bigr)^{\frac{1}{p}}
\]
turns this into approximation of indicators in $L^p$.
Since $f\mathbf1_{X_k}\to f$ in $L^p(X)$, density of simple functions
and approximation of their coefficients by rational numbers
---complex rational numbers in the complex case--- show that finite
sums of indicators $\mathbf1_{C\cap X_k}$ with these coefficients
form a countable dense set. This proves separability.
\end{proof}

\begin{corollary}[$L^p(M,\mathbf E)$ duality for Hermitian bundles]
\label{cor:dualidad-Lp-haz-hermitiano}
Let $(M,\mathbf{g})$ be a Riemannian manifold with or without boundary satisfying the second axiom of countability, and let $\mathbf{E}\to M$ be a
smooth complex Hermitian bundle of finite rank. For $1\leq p<\infty$, the map
\[
 L^{p'}(M,\overline{\mathbf{E}})\longrightarrow(L^p(M,\mathbf{E}))',
 \qquad
 \overline{\mathbf{v}}\longmapsto
 \left[u\longmapsto
 \int_M\mathcal R_{\mathbf{E}}(\overline{\mathbf{v}})(u)\,d\lambda_{\mathbf{g}}\right]
\]
is a complex-linear isometric isomorphism.
\end{corollary}

\begin{proof}
Write $\mu:=\lambda_{\mathbf g}$. Second countability, local
finiteness of $\mu$ and a cover by trivializations provide a
measurable partition $M=\bigsqcup_{j\geq1}A_j$, with a measurable
unitary frame $e_{j,1},\ldots,e_{j,r}$ on each piece.
Also choose increasing measurable sets $X_k$ of finite measure
that cover $M$. In each frame, if
$u=\sum_{\alpha=1}^r u_\alpha e_{j,\alpha}$, its fiber norm is
$(\sum_\alpha|u_\alpha|^2)^{\displaystyle\frac{1}{2}}$.

Let $\Lambda\in(L^p(M,\mathbf E))'$. For each $j,\alpha$, the map
that sends $f\in L^p(A_j)$ to $\Lambda(fe_{j,\alpha}\mathbf1_{A_j})$
is a continuous functional. Theorem~\ref{teo:dualidad-reflexividad-Lp}
gives a unique density $h_{j,\alpha}\in L^{p'}(A_j)$ representing
it as $f\mapsto\int_{A_j}fh_{j,\alpha}\,d\mu$.
On $A_j$, define
\[
 v_j:=\sum_{\alpha=1}^r\overline{h_{j,\alpha}}\,e_{j,\alpha}.
\]
There are finitely many measurable densities, so
$v_j\in L^{p'}(A_j,\mathbf E)$. For $u$ supported in $A_j$,
summing the coordinate representations gives
\[
 \Lambda(u)=\int_{A_j}\langle u,v_j\rangle_{\mathbf E}\,d\mu.
\]
Glue the $v_j$ into a measurable section $v$ and put
$q(x):=|v(x)|_{\mathbf E}$, finite almost everywhere.
The preceding representation holds for every section supported
in finitely many pieces $A_j$.

Suppose $1<p<\infty$. Let
\[
 C_k:=X_k\cap\left(\bigcup_{j=1}^k A_j\right)\cap\{q\leq k\},
 \qquad
 u_k:=\mathbf1_{C_k}q^{p'-2}v,
\]
with $u_k=0$ where $q=0$. Since $|u_k|_{\mathbf E}=q^{p'-1}$
on $C_k$, $u_k\in L^p$ and is supported in finitely many pieces.
Writing $I_k:=\int_{C_k}q^{p'}\,d\mu$, we obtain
\[
 I_k=\Lambda(u_k)\leq\|\Lambda\|\,\|u_k\|_{L^p}
     =\|\Lambda\|\,I_k^{\frac{1}{p}}.
\]
If $I_k>0$, division gives $I_k^{\displaystyle\frac{1}{p'}}\leq\|\Lambda\|$; if
$I_k=0$, the bound is immediate. The sets $C_k$ increase and
cover $M$ up to a null set. Monotone convergence proves
$v\in L^{p'}(M,\mathbf E)$ and $\|v\|_{L^{p'}}\leq\|\Lambda\|$.

If $p=1$, suppose $\{q>\|\Lambda\|+\varepsilon\}$ has positive
measure for some $\varepsilon>0$. The partition and exhaustion
give $j,k$ such that
\[
 A:=A_j\cap X_k\cap\{q>\|\Lambda\|+\varepsilon\},
 \qquad 0<\mu(A)<\infty.
\]
Define $u=\displaystyle\frac{v}{q}$ on $A$ and $u=0$ outside $A$. This section
belongs to $L^1$, is supported in $A_j$, and satisfies
\[
 \Lambda(u)=\int_A q\,d\mu
 \geq(\|\Lambda\|+\varepsilon)\mu(A)
 >\|\Lambda\|\,\|u\|_{L^1},
\]
a contradiction. Hence $q\leq\|\Lambda\|$ almost everywhere
and $v\in L^\infty$ with the same bound.

For either exponent, fiberwise Cauchy--Schwarz and Hölder give
\[
 \left|\int_M\langle u,v\rangle_{\mathbf E}\,d\mu\right|
 \leq\|u\|_{L^p}\|v\|_{L^{p'}}.
\]
The sections $u\mathbf1_{\bigcup_{j=1}^k A_j}$ converge to $u$
in $L^p$, since $p<\infty$. Passing to the limit in their
representations yields
$\Lambda(u)=\int_M\langle u,v\rangle_{\mathbf E}\,d\mu$
for every $u\in L^p$. The last bound and the one obtained for
$v$ imply $\|\Lambda\|=\|v\|_{L^{p'}}$.
Uniqueness of the scalar densities proves uniqueness of $v$.
Finally,
$\mathcal R_{\mathbf E}(\overline v)(u)=\langle u,v\rangle_{\mathbf E}$.
Complex linearity of $\mathcal R_{\mathbf E}$ on the conjugate
space, the representation and equality of norms prove the
asserted isomorphism.
\end{proof}
\begin{proposition}[Hölder's inequality] \label{desigualdad de holder}\index{Holder's inequality@Hölder's inequality}
Let $(X,\mathcal{A},\mu)$ be a measure space, and consider $r,p_{1},\dots,p_{m}\in [1,\infty]$ such that $\displaystyle\frac{1}{p_{1}}+\dots+\displaystyle\frac{1}{p_{m}}=\displaystyle\frac{1}{r}$, with the convention $\displaystyle\frac{1}{\infty}=0$. Then for any $f_{j}\in L^{p_{j}}(X)$, $1\leq j\leq m$, we have $\displaystyle\prod_{j=1}^{m}f_{j}\in L^{r}(X)$ and
\[
\left\|\displaystyle\prod_{j=1}^{m}f_{j}\right\|_{L^r(X)}
\leq \displaystyle\prod_{j=1}^{m}\|f_{j}\|_{L^{p_j}(X)}.
\]
\end{proposition}
\begin{proof}
Factors with $p_j=\infty$ are taken out using their norm in $L^\infty(X)$, so it suffices to consider indices with $p_j<\infty$. If any of them satisfies $\|f_j\|_{L^{p_j}(X)}=0$, the assertion is immediate. Otherwise, set $q_j:=\displaystyle\frac{p_j}{r}$; then $q_j\geq1$ and $\displaystyle\sum_{j=1}^{m}\displaystyle\frac{1}{q_j}=1$. If some $q_j=1$, it is the only finite index and the assertion is immediate. In the remaining case, the generalized numerical Young inequality,
\[
\displaystyle\prod_j a_j
\leq\displaystyle\sum_{j=1}^{m}\frac{a_j^{q_j}}{q_j},
\qquad a_j\geq0,
\]
follows by iterating $ab\leq\displaystyle\frac{a^q}{q}+\displaystyle\frac{b^{q'}}{q'}$. To prove the latter, with $b\geq0$ fixed, minimize $t\mapsto\displaystyle\frac{t^q}{q}-tb+\displaystyle\frac{b^{q'}}{q'}$; its minimum, attained when $t^{q-1}=b$, is zero. Applying the generalized inequality to
\[
a_j(x):=\frac{|f_j(x)|^r}{\|f_j\|_{L^{p_j}(X)}^r}
\]
and integrating gives
\[
\int_X\displaystyle\prod_j a_j\,d\mu
\leq\displaystyle\sum_{j=1}^{m}\frac{1}{q_j}\int_Xa_j^{q_j}\,d\mu
=\displaystyle\sum_{j=1}^{m}\frac{1}{q_j}=1.
\]
Multiplying by the norms taken out and taking the $r$th root proves the result. The case in which all exponents are infinite is the pointwise inequality for the essential supremum.
\end{proof}
\begin{proposition}\label{interpolacion 2}
Let $(X,\mathcal{A},\mu)$ be a finite measure space, and let $p\in(0,\infty)$. If $r\in(0,p)$ and $f\in\mathcal{L}^{p}(X)$, then $f\in \mathcal{L}^{r}(X)$ and
\[
 \|f\|_{\mathcal L^r(X)}
 \leq \|f\|_{\mathcal L^p(X)}\mu(X)^{\frac{1}{r}-\frac{1}{p}}.
 \]
\end{proposition}
\begin{proof}
Apply Proposition~\ref{desigualdad de holder} to $|f|^r$ and the constant function $1$, with exponents $\displaystyle\frac{p}{r}$ and $\displaystyle\frac{p}{p-r}$. Thus
\[
\int_X|f|^r\,d\mu
\leq
\left(\int_X|f|^p\,d\mu\right)^{\frac{r}{p}}
\mu(X)^{1-\frac{r}{p}}.
\]
Taking the $r$th root gives the stated inequality.
\end{proof}
\begin{proposition}[Jensen's inequality]\label{desigualdad de Jensen}\index{Jensen's inequality} Let $(X,\mathcal{A},\mu)$ be a finite measure space with $\mu(X)=1$, let $\phi\colon \mathbb{R}\longrightarrow\mathbb{R}$ be convex, and let $f\in \mathcal{L}^{1}(X)$. Then, whenever the right-hand side is well defined, \[\phi\left(\int_{X}fd\mu\right)\leq \int_{X}\phi(f)d\mu\]
\end{proposition}
\begin{proof}
Set $a:=\displaystyle\int_Xf\,d\mu$. Convexity of $\phi$ implies that there exists a supporting slope $c\in\mathbb R$ at $a$: indeed, any number between the one-sided slopes
\[
\sup_{x\in(-\infty,a)}\frac{\phi(a)-\phi(x)}{a-x}
\quad\text{and}\quad
\inf_{y>a}\frac{\phi(y)-\phi(a)}{y-a}
\]
satisfies $\phi(t)\geq\phi(a)+c(t-a)$ for every $t\in\mathbb R$. Applying this inequality to $t=f(x)$ and integrating gives
\[
\int_X\phi(f)\,d\mu
\geq\phi(a)+c\left(\int_Xf\,d\mu-a\mu(X)\right)
=\phi(a).
\]
\end{proof}
\begin{proposition}\label{ejericio 14.66}
Let $\omega\subseteq \Omega\subseteq \mathbb{R}^{n}$ be open sets, and let $p\in[1,\infty]$. Then:
\begin{enumerate}[label=(\alph*)]
\item If $f\in L^{p}(\Omega)$, then $f\mathbf{1}_{\omega}\in L^{p}(\omega)$ and \[\|f\mathbf{1}_{\omega}\|_{L^{p}(\omega)}\leq \|f\|_{L^{p}(\Omega)}\]
\item If $f\in L^{p}(\Omega)$ and $X$ is a measurable subset of $\Omega$ with $\lambda_{n}(X)<\infty$, then:
\begin{enumerate}
\item if $p<\infty$, we have $f\mathbf{1}_{X}\in L^{r}(\Omega)$ for every $r\in [1,p]$ and \[\|f\mathbf{1}_{X}\|_{L^{r}(\Omega)}\leq \lambda_{n}(X)^{\frac{p-r}{rp}}\|f\|_{L^{p}(\Omega)};\]
\item if $p=\infty$, we have $f\mathbf{1}_{X}\in L^{r}(\Omega)$ for every $r\in[1,\infty)$ and \[\|f\mathbf{1}_{X}\|_{L^{r}(\Omega)}\leq \lambda_n(X)^{\frac{1}{r}}\|f\|_{L^\infty(\Omega)}.\]
\end{enumerate}
\item If $f_{k}\to f$ in $L^{p}(\Omega)$ and $\lambda_n(\omega)<\infty$, then:
\begin{enumerate}
\item if $p<\infty$, $f_{k}\mathbf{1}_{\omega}\to f\mathbf{1}_{\omega}$ in $L^{r}(\omega)$ for every $r\in [1,p]$;
\item if $p=\infty$, $f_{k}\mathbf{1}_{\omega}\to f\mathbf{1}_{\omega}$ in $L^{r}(\omega)$ for every $r\in[1,\infty)$.
\end{enumerate}
\end{enumerate}
\end{proposition}
\begin{proof}
Part (a) follows directly by integrating over $\omega\subseteq\Omega$. For (b), if $1\leq r<p<\infty$, apply Proposition~\ref{desigualdad de holder} on $X$ to $|f|^r$ and $1$, with exponents $\displaystyle\frac{p}{r}$ and $\displaystyle\frac{p}{p-r}$; if $r=p$, the assertion is immediate. For $p=\infty$, use $|f|\leq\|f\|_{L^\infty(\Omega)}$ almost everywhere on $X$. This gives exactly the two bounds in (b). Finally, apply these bounds to $(f_k-f)\mathbf1_\omega$ to obtain (c).
\end{proof}
\begin{proposition}\label{ejercicio 14.70}
If $f=(f_{1},\dots,f_{m})$ and $f_{1},\dots,f_{m}\in L^{s}(\Omega)$, define
\[
 \|f\|_{(L^s(\Omega))^m}
 :=\begin{cases}
 \left(\displaystyle\sum_{j=1}^{m}\|f_j\|_{L^s(\Omega)}^s\right)^{\frac{1}{s}},&1\leq s<\infty,\\
 \displaystyle\max_{1\leq j\leq m}\|f_j\|_{L^\infty(\Omega)},&s=\infty.
 \end{cases}
 \]
If $\lambda_{n}(\Omega)<\infty$, $1\leq p<s\leq \infty$, and $f_{1},\dots,f_{m}\in L^{s}(\Omega)$, then
\[\begin{cases}
 \|f\|_{(L^p(\Omega))^m}\leq (m\lambda_{n}(\Omega))^{\frac{s-p}{sp}}\|f\|_{(L^s(\Omega))^m} & \text{if $s\in[1,\infty)$},\\
 \|f\|_{(L^p(\Omega))^m}\leq (m\lambda_{n}(\Omega))^{\frac{1}{p}}\|f\|_{(L^\infty(\Omega))^m} & \text{if $s=\infty$}.
 \end{cases}\]
\end{proposition}
\begin{proof}
Consider $\Omega\times\{1,\dots,m\}$ with the product of Lebesgue measure and counting measure, and define $F(x,j):=f_j(x)$. This space has measure $m\lambda_n(\Omega)$, and
\[
\|F\|_{L^t(\Omega\times\{1,\dots,m\})}
=\|f\|_{(L^t(\Omega))^m}
\]
for $1\leq t\leq\infty$. If $s<\infty$, Proposition~\ref{desigualdad de holder}, applied to $|F|^p$ and $1$ with exponents $\displaystyle\frac{s}{p}$ and $\displaystyle\frac{s}{s-p}$, gives the first inequality. The second follows by integrating the bound $|F|\leq\|F\|_{L^\infty(\Omega\times\{1,\dots,m\})}$.
\end{proof}
\begin{proposition}\label{ejercicio 14.69}
If $0<\lambda_{n}(\Omega)<\infty$ and $f\in L^{\infty}(\Omega)$, then
\[
 \|f\|_{L^\infty(\Omega)}
 =\lim_{p\to\infty}\lambda_{n}(\Omega)^{-\frac{1}{p}}
 \|f\|_{L^p(\Omega)}.
 \]
\end{proposition}
\begin{proof}
The inequality
\[
\lambda_n(\Omega)^{-\frac{1}{p}}\|f\|_{L^p(\Omega)}
\leq\|f\|_{L^\infty(\Omega)}
\]
gives the upper bound. If $0<a<\|f\|_{L^\infty(\Omega)}$, the set $A_a:=\{x\in\Omega\mid |f(x)|>a\}$ has positive measure and
\[
\lambda_n(\Omega)^{-\frac{1}{p}}\|f\|_{L^p(\Omega)}
\geq a\left(\frac{\lambda_n(A_a)}{\lambda_n(\Omega)}\right)^{\frac{1}{p}}.
\]
The right-hand side converges to $a$ as $p\to\infty$. Then letting $a\to\bigl(\|f\|_{L^\infty(\Omega)}\bigr)^{-}$ gives the lower bound and hence the limit.
\end{proof}
\begin{proposition}\label{Teorema 14.28}
Let $(f_{k})$ be a sequence in $L^{p}(\Omega)$ such that $f_{k}\to f$ in $L^{p}(\Omega)$.
\begin{enumerate}[label=(\alph*)]
\item If $p\in[1,\infty)$, there exist a subsequence $(f_{k_{j}})$ of $(f_{k})$ and a function $g\in L^{p}(\Omega)$ such that $f_{k_{j}}(x)\to f(x)$ for almost every $x\in\Omega$, and $|f_{k_{j}}|\leq g(x)$ for almost every $x\in\Omega$, for every $j\in\mathbb{N}$.
\item If $p=\infty$, then $f_{k}(x)\to f(x)$ for almost every $x\in\Omega$, and there exists $c\in\mathbb{R}$ such that $|f_{k}(x)|\leq c$ for almost every $x\in\Omega$ and every $k\in\mathbb{N}$.
\end{enumerate}
\end{proposition}
\begin{proof}
First suppose $p<\infty$. Choosing a subsequence inductively, we may ensure that
\[
\|f_{k_{j+1}}-f_{k_j}\|_{L^p(\Omega)}\leq2^{-j}
\qquad(j\in\mathbb N).
\]
Minkowski's inequality and Theorem~\ref{convergencia monotona} show that
\[
g:=|f_{k_1}|+\displaystyle\sum_{j\in\mathbb N}|f_{k_{j+1}}-f_{k_j}|
\]
belongs to $L^p(\Omega)$. The series of differences converges absolutely almost everywhere; hence $(f_{k_j})$ converges almost everywhere to a function $h$, and $|f_{k_j}|\leq g$. Moreover, $|h|\leq g$ almost everywhere. The dominated convergence theorem~\ref{convergencia dominada}, applied to $|f_{k_j}-h|^p\leq(2g)^p$, shows that $f_{k_j}\to h$ in $L^p(\Omega)$. Since also $f_{k_j}\to f$ in $L^p(\Omega)$, uniqueness of the limit gives $h=f$ almost everywhere.

If $p=\infty$, choose representatives and a common null set outside which
\[
|f_k(x)-f(x)|\leq\|f_k-f\|_{L^\infty(\Omega)}\longrightarrow0.
\]
The sequence of norms $\|f_k-f\|_{L^\infty(\Omega)}$ is bounded; this gives a constant $c$ dominating all the $|f_k|$ outside a null set.
\end{proof}

\begin{theorem}[Interpolation inequality]\label{Ejercicio 14.67}\index{interpolation inequality}
Let $1\leq p<s<r\leq\infty$. If $f\in L^{p}(\Omega)\cap L^{r}(\Omega)$, then $f\in L^{s}(\Omega)$ and
\[
 \|f\|_{L^s(\Omega)}
 \leq \|f\|_{L^p(\Omega)}^{1-\alpha}
 \|f\|_{L^r(\Omega)}^{\alpha},
 \]
where $\alpha\in (0,1)$ satisfies $\displaystyle\frac{1}{s}=\displaystyle\frac{1-\alpha}{p}+\displaystyle\frac{\alpha}{r}$ if $r<\infty$ and $\alpha:=1-\displaystyle\frac{p}{s}$ if $r=\infty$.
\end{theorem}
\begin{proof}
If $r<\infty$, apply Proposition~\ref{desigualdad de holder} to
\[
|f|^s=|f|^{(1-\alpha)s}|f|^{\alpha s}
\]
with exponents $\displaystyle\frac{p}{(1-\alpha)s}$ and $\displaystyle\frac{r}{\alpha s}$, which are conjugate by the relation between $p,s,r$ and $\alpha$. Taking the $s$th root yields the inequality. If $r=\infty$, then $s=(1-\alpha)^{-1}p$ and
\[
\int_\Omega|f|^s\,d\lambda_n
\leq\|f\|_{L^\infty(\Omega)}^{s-p}
\int_\Omega|f|^p\,d\lambda_n,
\]
which gives the same conclusion.
\end{proof}
\begin{proposition}\label{Proposición 14.31}
If $\lambda_{n}(\Omega)<\infty$ and $1\leq p<s\leq \infty$, then $L^{s}(\Omega)\subseteq L^{p}(\Omega)$ and this inclusion is continuous. Moreover, for every $f\in L^{s}(\Omega)$, \[\begin{cases}
 \|f\|_{L^p(\Omega)}\leq (\lambda_{n}(\Omega))^{\frac{s-p}{sp}}\|f\|_{L^s(\Omega)} & \text{if $s\in[1,\infty)$},\\
 \|f\|_{L^p(\Omega)}\leq (\lambda_{n}(\Omega))^{\frac{1}{p}}\|f\|_{L^\infty(\Omega)} & \text{if $s=\infty$}.
 \end{cases}\]
\end{proposition}
\begin{proof}
If $s<\infty$, apply Proposition~\ref{desigualdad de holder} to $|f|^p$ and $1$ with exponents $\displaystyle\frac{s}{p}$ and $\displaystyle\frac{s}{s-p}$. If $s=\infty$, integrate the inequality $|f|^p\leq\|f\|_{L^\infty(\Omega)}^p$. The resulting bounds also prove continuity of the inclusion.
\end{proof}
\begin{proposition}\label{densidad de Cc en Lp}
$C_{c}^{\infty}(\Omega)$ is dense in $L^{p}(\Omega)$ for every $p\in [1,\infty)$.
\end{proposition}
\begin{proof}
By Proposition~\ref{prop:densidad-funciones-simples-Lp}, it suffices to approximate in $L^p(\Omega)$ the indicator of a measurable set $A\subseteq\Omega$ of finite measure. Regularity of Lebesgue measure in Theorem~\ref{teorema 5.5} allows us to choose a compact set $K\subseteq A$ and an open set $U$ with $K\subseteq U\Subset\Omega$ such that
\[
\lambda_n(A\setminus K)+\lambda_n(U\setminus K)<\varepsilon^p.
\]
By Proposition~\ref{funciones flan}, there exists $\eta\in C_c^\infty(\Omega)$ such that $0\leq\eta\leq1$, $\eta=1$ on $K$, and $\operatorname{supp}(\eta)\subseteq U$. Then
\[
\|\mathbf1_A-\eta\|_{L^p(\Omega)}^p
\leq\lambda_n(A\setminus K)+\lambda_n(U\setminus K)
<\varepsilon^p.
\]
Approximating each term of a simple function and using the triangle inequality proves density.
\end{proof}
\begin{proposition}\label{corolario 14.47 relativamente compacto en Lp}
Let $p\in [1,\infty)$, let $\Omega\subseteq \mathbb{R}^{n}$ be open, and let $\mathcal{K}\subseteq L^{p}(\Omega)$ be bounded with the following properties:
\begin{enumerate}[label=(\alph*)]
\item For each $\varepsilon>0$ and each open set $\omega\Subset \Omega$, there exists $\delta\in (0,d(\omega,\mathbb{R}^{n}\setminus \Omega))$ such that $\|T_{\boldsymbol{\xi}}f-f\|_{L^{p}(\omega)}<\varepsilon$ for every $\boldsymbol{\xi}\in\mathbb{R}^{n}$ with $\|\boldsymbol{\xi}\|<\delta$ and every $f\in\mathcal{K}$.
\item For each $\varepsilon>0$, there exists an open set $\omega\Subset\Omega$ such that $\|f\|_{L^{p}(\Omega\setminus \omega)}<\varepsilon$ for every $f\in\mathcal{K}$.
\end{enumerate}
Here $T_{\boldsymbol{\xi}}f(x):=f(x+\boldsymbol{\xi})$ whenever
$x,x+\boldsymbol{\xi}\in\Omega$.
Then $\mathcal{K}$ is relatively compact in $L^{p}(\Omega)$.
\end{proposition}
\begin{proof}
Fix $\varepsilon>0$. By (b), choose $\omega\Subset\Omega$ so that
$\|f\|_{L^p(\Omega\setminus\omega)}<\varepsilon$ for every $f\in\mathcal K$, and choose $\omega_1$ with $\omega\Subset\omega_1\Subset\Omega$. For small $h>0$, take a finite family of disjoint cubes $(Q_\ell)_{\ell=1}^N$, of side length $h$, contained in $\omega_1$, whose union $U_h$ contains $\omega$. Define
\[
P_hf
:=
\displaystyle\sum_{\ell=1}^N
\left(\frac{1}{|Q_\ell|}\int_{Q_\ell}f(y)\,dy\right)\mathbf1_{Q_\ell}.
\]
Jensen's inequality and Theorem~\ref{teo:tonelli} give
\[
\|P_hf-f\|_{L^p(U_h)}^p
\leq
\frac{1}{h^n}\displaystyle\sum_{\ell=1}^N
\int_{Q_\ell}\int_{Q_\ell}|f(y)-f(x)|^p\,dy\,dx.
\]
Writing $y=x+z$ and observing that $\|z\|\leq\sqrt n h$, the right-hand side is bounded by
\[
C h^{-n}\int_{\{\|z\|\leq\sqrt n h\}}
\|T_zf-f\|_{L^p(\omega_1)}^p\,dz.
\]
Condition (a) makes this quantity less than $C\varepsilon^p$ uniformly in $f\in\mathcal K$ if $h$ is sufficiently small. Since $P_hf=0$ outside $U_h$ and $\Omega\setminus U_h\subseteq\Omega\setminus\omega$, we obtain
$\|P_hf-f\|_{L^p(\Omega)}\leq C\varepsilon$.

The set $P_h(\mathcal K)$ is bounded in the finite-dimensional space spanned by $(\mathbf1_{Q_\ell})_{\ell=1}^N$ and is therefore totally bounded. The preceding uniform approximation shows that $\mathcal K$ is totally bounded in the complete space $L^p(\Omega)$; its closure is compact.
\end{proof}

\section{Integral inequalities and tools from harmonic analysis}
\label{ap:herramientas-analisis-armonico}

In this section we collect several tools that will appear repeatedly in the Littlewood--Paley decomposition and the study of Besov and Triebel--Lizorkin spaces. The first part contains integral and discrete inequalities allowing norms, integrals, and sums over scales to be interchanged. We then introduce Rademacher functions and the Hardy--Littlewood and Fefferman--Stein maximal estimates. Finally, we record the Calderón--Zygmund decomposition and Peetre maximal control, which allow families of frequency-localized functions to be treated uniformly.

Throughout this section, $\mathcal M$ denotes the uncentered Hardy--Littlewood maximal operator on $\mathbb R^n$.

\begin{proposition}[Cavalieri and Chebyshev]\label{prop:cavalieri-chebyshev}
If $f\colon X\longrightarrow[0,\infty]$ is measurable and $0<p<\infty$, then
\[
\int_Xf^p\,d\mu
=
p\int_0^\infty\lambda^{p-1}\mu(\{f>\lambda\})\,d\lambda.
\]
In particular, if $f\in L^p(X)$ and $\lambda>0$, then
\[
\mu(\{|f|>\lambda\})
\leq
\lambda^{-p}\|f\|_{L^p(X)}^p.
\]
\end{proposition}

\begin{proof}
First suppose that $f=\sum_{i=1}^Na_i\mathbf1_{A_i}$ is simple,
with $a_i>0$ and disjoint measurable sets $A_i$. Integrating a finite
sum of nonnegative functions gives
\[
\begin{aligned}
p\int_0^\infty\lambda^{p-1}\mu(\{f>\lambda\})\,d\lambda
&=\sum_{i=1}^N\mu(A_i)
       \int_0^{a_i}p\lambda^{p-1}\,d\lambda\\
&=\sum_{i=1}^Na_i^p\mu(A_i)=\int_Xf^p\,d\mu.
\end{aligned}
\]
The identity includes infinite measures of the $A_i$.
For general $f$, choose nonnegative simple functions $f_N\uparrow f$.
Then $f_N^p\uparrow f^p$ and, for each $\lambda>0$,
$\{f_N>\lambda\}\uparrow\{f>\lambda\}$. Continuity from below
of the measure and Theorem~\ref{convergencia monotona}, applied on
$X$ and on $(0,\infty)$, allow passage to the limit in the identity
for $f_N$. This proves Cavalieri's formula on any measure space,
without a $\sigma$-finiteness assumption.

For Chebyshev's inequality,
\[
\|f\|_{L^p(X)}^p
\geq\int_{\{|f|>\lambda\}}|f|^p\,d\mu
\geq\lambda^p\mu(\{|f|>\lambda\}).
\]
Division by $\lambda^p>0$ gives the asserted bound.
\end{proof}

\begin{theorem}[Minkowski's integral inequality]\label{teo:minkowski-integral}
Let $(X,\mu)$ and $(Y,\nu)$ be $\sigma$-finite measure spaces, and let $1\leq p\leq\infty$. If $F\colon X\times Y\longrightarrow\mathbb C$ is measurable and
\[
\int_Y\|F(\,\cdot,y)\|_{L^p(X)}\,d\nu(y)<\infty,
\]
then $\displaystyle x\mapsto\int_YF(x,y)\,d\nu(y)$ belongs to $L^p(X)$ and
\[
\left\|\int_YF(\,\cdot,y)\,d\nu(y)\right\|_{L^p(X)}
\leq
\int_Y\|F(\,\cdot,y)\|_{L^p(X)}\,d\nu(y).
\]
\end{theorem}

\begin{proof}
Put $m(y):=\|F(\,\cdot,y)\|_{L^p(X)}$ and
$A:=\int_Ym(y)\,d\nu(y)<\infty$. This function is measurable:
for $p<\infty$, this follows from Tonelli applied to $|F|^p$;
for $p=\infty$, for every $t\geq0$ we have
\[
 \{y\mid m(y)>t\}
 =\left\{y\middle|\mu\{x\mid |F(x,y)|>t\}>0\right\},
\]
and Tonelli applied to indicators proves measurability.
Define the nonnegative measurable function, initially with values
in $[0,\infty]$,
\[
 H(x):=\int_Y|F(x,y)|\,d\nu(y).
\]
For $p=1$, Tonelli gives $\int_XH\,d\mu=A$.

Suppose $1<p<\infty$. Choose increasing measurable sets
$X_k$ of finite measure covering $X$, and put
$H_k:=\min\{H,k\}\mathbf1_{X_k}$. Then $H_k\in L^p(X)$
and $H_k\uparrow H$. If $\|H_k\|_{L^p(X)}>0$, define
\[
 g_k:=\frac{H_k^{p-1}}{\|H_k\|_{L^p(X)}^{p-1}},
 \qquad \|g_k\|_{L^{p'}(X)}=1.
\]
Tonelli and Hölder, applied to each section, give
\[
 \begin{aligned}
 \|H_k\|_{L^p(X)}
 &=\int_Xg_kH_k\,d\mu
 \leq\int_Xg_kH\,d\mu\\
 &=\int_Y\int_Xg_k(x)|F(x,y)|\,d\mu(x)\,d\nu(y)
 \leq\int_Ym(y)\,d\nu(y)=A.
 \end{aligned}
\]
If $\|H_k\|_{L^p(X)}=0$, the same bound is immediate.
Monotone convergence applied to $H_k^p$ proves that
$H\in L^p(X)$ and $\|H\|_{L^p(X)}\leq A$.

If $p=\infty$, the measurable set
$N:=\{(x,y)\mid |F(x,y)|>m(y)\}$ has $\mu$-null horizontal
sections. Tonelli implies that its vertical sections are
$\nu$-null for $\mu$-almost every $x$. For these $x$,
\[
 H(x)\leq\int_Ym(y)\,d\nu(y)=A.
\]
Thus $H\in L^\infty(X)$ and its norm is at most $A$.

In every case, $H$ is finite almost everywhere. Hence the integral
$h(x):=\int_YF(x,y)\,d\nu(y)$ converges absolutely for almost
every $x$. Define it to be zero on the measurable exceptional set.
Integrating the positive and negative parts of the real and imaginary
components shows that $h$ is measurable.
Since $|h|\leq H$, we obtain $h\in L^p(X)$ and the stated inequality.
\end{proof}

\begin{proposition}[Schur's test]
\label{prop:criterio-schur-operadores-integrales}
\index{Schur's test}
Let $(X,\mu)$ and $(Y,\nu)$ be $\sigma$-finite measure spaces, let
$r,s\in\mathbb N$, and let
\[
K\colon X\times Y
\longrightarrow
\mathcal L(\mathbb C^r,\mathbb C^s)
\]
be a measurable kernel. Suppose that there exist $A,B\geq0$ such that
\[
\sup_{x\in X}
\int_Y
\|K(x,y)\|_{\mathcal L(\mathbb C^r,\mathbb C^s)}\,d\nu(y)
\leq A
\]
and
\[
\sup_{y\in Y}
\int_X
\|K(x,y)\|_{\mathcal L(\mathbb C^r,\mathbb C^s)}\,d\mu(x)
\leq B.
\]
Then, for each $f\in L^2(Y,\mathbb C^r)$, the integral
\[
(T_Kf)(x):=\int_YK(x,y)f(y)\,d\nu(y)
\]
exists for almost every $x\in X$ and defines a continuous linear operator
\[
T_K\colon L^2(Y,\mathbb C^r)
\longrightarrow L^2(X,\mathbb C^s)
\]
satisfying
\[
\|T_K\|_{\mathcal L(
L^2(Y,\mathbb C^r),L^2(X,\mathbb C^s))}
\leq\sqrt{AB}.
\]
\end{proposition}

\begin{proof}
Set
$k(x,y):=\|K(x,y)\|_{\mathcal L(\mathbb C^r,\mathbb C^s)}$.
For $f\in L^2(Y,\mathbb C^r)$, Cauchy--Schwarz with respect to the measure
$k(x,y)\,d\nu(y)$ gives
\begin{equation}
\left(\int_Yk(x,y)\|f(y)\|_{\mathbb C^r}\,d\nu(y)\right)^2
\leq
A\int_Yk(x,y)\|f(y)\|_{\mathbb C^r}^2\,d\nu(y).
\label{eq:schur-existencia-integral}
\end{equation}
Tonelli's theorem~\ref{teo:tonelli} and the second hypothesis show that
\[
\begin{aligned}
\int_X\int_Y
k(x,y)\|f(y)\|_{\mathbb C^r}^2\,d\nu(y)\,d\mu(x)
&=
\int_Y\|f(y)\|_{\mathbb C^r}^2
\left(\int_Xk(x,y)\,d\mu(x)\right)d\nu(y)\\
&\leq B\|f\|_{L^2(Y,\mathbb C^r)}^2.
\end{aligned}
\]
Thus the last factor in
\eqref{eq:schur-existencia-integral} is finite for almost every $x$; the
integral defining $T_Kf(x)$ is absolutely convergent in
$\mathbb C^s$ at these points. Measurability of $T_Kf$ follows
by integrating measurable functions componentwise. Moreover,
\[
\begin{aligned}
\|T_Kf\|_{L^2(X,\mathbb C^s)}^2
&\leq
A\int_X\int_Y
k(x,y)\|f(y)\|_{\mathbb C^r}^2\,d\nu(y)\,d\mu(x)\\
&\leq AB\|f\|_{L^2(Y,\mathbb C^r)}^2.
\end{aligned}
\]
This proves the estimate and, in particular, continuity. The same
proof applies to a measurable kernel
$K(x,y)\colon E_y\longrightarrow F_x$ between measurable Hermitian bundles of finite
rank, replacing Euclidean norms by bundle norms.
\end{proof}

\begin{proposition}[Young's inequality for sequences]\label{prop:young-sucesiones}
Let $1\leq q\leq\infty$, $a=(a_j)_{j\in\mathbb Z}\in\ell^q(\mathbb Z)$, and $b=(b_j)_{j\in\mathbb Z}\in\ell^1(\mathbb Z)$. If
\[
(a*b)_k:=\sum_{j\in\mathbb Z}a_jb_{k-j},
\qquad k\in\mathbb Z,
\]
then
\[
\|a*b\|_{\ell^q(\mathbb Z)}
\leq
\|a\|_{\ell^q(\mathbb Z)}\|b\|_{\ell^1(\mathbb Z)}.
\]
\end{proposition}

\begin{proof}
Every sequence in $\ell^q(\mathbb Z)$ is bounded, and
$\|a\|_{\ell^\infty}\leq\|a\|_{\ell^q}$. Hence, for each $k$,
\[
\sum_{j\in\mathbb Z}|a_jb_{k-j}|
\leq\|a\|_{\ell^\infty}\|b\|_{\ell^1}<\infty.
\]
This justifies the definition of convolution and proves the bound
when $q=\infty$. If $q=1$, Tonelli for counting measure gives
\[
\sum_{k\in\mathbb Z}|(a*b)_k|
\leq\sum_{k\in\mathbb Z}\sum_{j\in\mathbb Z}|a_j||b_{k-j}|
=\|a\|_{\ell^1}\|b\|_{\ell^1}.
\]
If $1<q<\infty$ and $b\neq0$, apply Hölder to
$|a_j||b_{k-j}|^{\displaystyle\frac{1}{q}}$ and
$|b_{k-j}|^{1-\displaystyle\frac{1}{q}}$, first in finite sums and then
by monotone convergence. We obtain
\[
|(a*b)_k|^q
\leq\|b\|_{\ell^1}^{q-1}
          \sum_{j\in\mathbb Z}|a_j|^q|b_{k-j}|.
\]
Summing in $k$ and using Tonelli again yields
\[
\|a*b\|_{\ell^q}^q\leq\|b\|_{\ell^1}^q\|a\|_{\ell^q}^q.
\]
If $b=0$, all assertions are immediate.
\end{proof}

\begin{theorem}[Discrete Hardy inequalities]\label{teo:desigualdades-hardy}
Let $\alpha>0$, $1\leq q\leq\infty$, and let $(a_j)_{j\in\mathbb Z}$ be a nonnegative sequence. Then
\[
\left\|\left(\sum_{\substack{j\in\mathbb Z\\j\leq k}}2^{-\alpha(k-j)}a_j\right)_{k\in\mathbb Z}\right\|_{\ell^q(\mathbb Z)}
+
\left\|\left(\sum_{\substack{j\in\mathbb Z\\j\geq k}}2^{-\alpha(j-k)}a_j\right)_{k\in\mathbb Z}\right\|_{\ell^q(\mathbb Z)}
\leq
C_{\alpha,q}\|(a_j)_{j\in\mathbb Z}\|_{\ell^q(\mathbb Z)}.
\]
\end{theorem}

\begin{proof}
The first sequence is the convolution of $(a_j)_{j\in\mathbb Z}$ with
\[
b_m:=2^{-\alpha m}\mathbf1_{\{m\geq0\}},
\qquad m\in\mathbb Z,
\]
and the second with $\widetilde b_m:=2^{\alpha m}\mathbf1_{\{m\leq0\}}$. Both sequences belong to $\ell^1(\mathbb Z)$. Proposition~\ref{prop:young-sucesiones} gives the two bounds.
\end{proof}

\begin{theorem}[Hardy's integral inequalities]
\label{teo:desigualdades-hardy-integrales}
Let $\alpha>0$, $1\leq q\leq\infty$, and let $f\colon (0,\infty)\longrightarrow[0,\infty]$ be measurable. Then
\[
\left\|
r^{-\alpha}\int_0^r t^\alpha f(t)\,\frac{dt}{t}
\right\|_{L^q((0,\infty),\frac{dr}{r})}
\leq
\frac{1}{\alpha}
\|f\|_{L^q((0,\infty),\frac{dt}{t})}
\]
and
\[
\left\|
r^\alpha\int_r^\infty t^{-\alpha}f(t)\,\frac{dt}{t}
\right\|_{L^q((0,\infty),\frac{dr}{r})}
\leq
\frac{1}{\alpha}
\|f\|_{L^q((0,\infty),\frac{dt}{t})}.
\]
The same assertions hold for Banach-valued functions upon replacing $f(t)$ by its norm.
\end{theorem}

\begin{proof}
Make the substitutions $r=e^x$, $t=e^y$, and set $g(y):=f(e^y)$. Since $\displaystyle\frac{dt}{t}=dy$ and $\displaystyle\frac{dr}{r}=dx$, the first expression becomes
\[
\int_{-\infty}^x e^{-\alpha(x-y)}g(y)\,dy
=
\bigl(e^{-\alpha\,\cdot}\mathbf1_{[0,\infty)}*g\bigr)(x),
\]
while the second is the convolution of $g$ with $e^{\alpha\,\cdot}\mathbf1_{(-\infty,0]}$. Both kernels have $L^1(\mathbb R)$ norm equal to $\displaystyle\frac{1}{\alpha}$. Young's inequality in Theorem~\ref{teo:young-convoluciones} completes the proof, including the endpoints $q=1$ and $q=\infty$.
\end{proof}

Rademacher functions allow a Euclidean norm of coefficients to be converted into a scalar integral. This randomization is the mechanism we shall later use to pass from multipliers with arbitrary signs to the Littlewood--Paley square function.

\begin{lemma}[Paley--Zygmund inequality]
\label{lem:paley-zygmund}
Let $(A,\mathcal A,\mu)$ be a probability space, let $X\in L^2(A)$ be nonnegative, and let $0<\vartheta<1$. Then
\[
\mu\{X>\vartheta\,\mathbb EX\}
\geq
(1-\vartheta)^2
\frac{(\mathbb EX)^2}{\mathbb E(X^2)}.
\]
Here $\displaystyle \mathbb EX:=\int_AX\,d\mu$; if $\mathbb E(X^2)=0$, both sides are interpreted as zero.
\end{lemma}

\begin{proof}
The case $\mathbb EX=0$ is immediate. Suppose $\mathbb EX>0$, and let $E:=\{X>\vartheta\mathbb EX\}$. Then
\[
\mathbb EX
=
\int_{A\setminus E}X\,d\mu+\int_EX\,d\mu
\leq
\vartheta\mathbb EX+
(\mathbb E(X^2))^{\frac{1}{2}}\mu(E)^{\frac{1}{2}},
\]
where we used Cauchy--Schwarz, Theorem~\ref{teo:cauchy-schwarz-hilbert}, in the last term. Subtracting $\vartheta\mathbb EX$, squaring, and dividing by $\mathbb E(X^2)$ gives the assertion.
\end{proof}

\begin{theorem}[Khintchine's inequality]\label{teo:khintchine}
Let $(r_j)_{j\in\mathbb N_0}$ be the Rademacher functions on $[0,1]$. For each $0<p<\infty$ there exist constants $A_p,B_p>0$ such that, for every finite set $J\subseteq\mathbb N_0$ and every family of scalars $(a_j)_{j\in J}$,
\[
A_p\left(\sum_{j\in J}|a_j|^2\right)^{\frac{1}{2}}
\leq
\left(\int_0^1\left|\sum_{j\in J}a_jr_j(t)\right|^pdt\right)^{\frac{1}{p}}
\leq
B_p\left(\sum_{j\in J}|a_j|^2\right)^{\frac{1}{2}}.
\]
\end{theorem}

\begin{proof}
First suppose that the coefficients are real, and set
\[
S:=\sum_{j\in J}a_jr_j,
\qquad
\sigma^2:=\sum_{j\in J}|a_j|^2.
\]
Independence of the Rademacher functions and the inequality $\cosh x\leq e^{\displaystyle\frac{x^2}{2}}$ imply
\[
\int_0^1e^{\lambda S(t)}dt
=
\prod_{j\in J}\cosh(\lambda a_j)
\leq
e^{\frac{\lambda^2\sigma^2}{2}}.
\]
Apply Chebyshev, Proposition~\ref{prop:cavalieri-chebyshev}, to $e^{\lambda S}$ and $e^{-\lambda S}$. For $u>0$ and $\lambda>0$, we obtain
\[
|\{t\mid S(t)>u\}|
\leq
e^{-\lambda u+\frac{\lambda^2\sigma^2}{2}},
\]
and the same bound for $-S$. If $\sigma>0$, choosing $\lambda=\displaystyle\frac{u}{\sigma^2}$ gives
\[
|\{t\in[0,1]\mid |S(t)|>u\}|
\leq
2e^{-\frac{u^2}{2\sigma^2}}.
\]
Cavalieri's formula then yields
\[
\int_0^1|S|^pdt
=
p\int_0^\infty u^{p-1}|\{|S|>u\}|\,du
\leq
C_p\sigma^p,
\]
which is the upper bound.

For the lower bound, by homogeneity we may assume $\sigma=1$. Independence gives
\[
\mathbb E(S^2)=1,
\qquad
\mathbb E(S^4)
=
\sum_{j\in J}a_j^4
+6\sum_{\substack{j,k\in J\\j<k}}a_j^2a_k^2
=
3-2\sum_{j\in J}a_j^4
\leq3.
\]
Apply Lemma~\ref{lem:paley-zygmund} to $X=S^2$ with $\vartheta=\displaystyle\frac{1}{2}$. Thus
\[
|\{|S|>2^{-\frac{1}{2}}\}|
\geq
\frac{1}{12},
\]
and consequently,
\[
\int_0^1|S|^pdt
\geq
\frac{1}{12}\,2^{-\frac{p}{2}}.
\]
This proves the lower bound for real coefficients. In the complex case, apply the real bounds to the real and imaginary parts. Both have sum of squares at most $\sigma^2$, and at least one has sum of squares no less than $\displaystyle\frac{\sigma^2}{2}$. The inequalities
\[
|\operatorname{Re}S|,|\operatorname{Im}S|
\leq|S|
\leq
|\operatorname{Re}S|+|\operatorname{Im}S|
\]
complete the proof.
\end{proof}

Before studying the maximal operator, we record the interpolation principle used to pass from weak endpoint estimates to strong estimates at intermediate exponents.

\begin{theorem}[Marcinkiewicz interpolation]\label{teo:interpolacion-marcinkiewicz}
Let $(X,\mu_X)$ and $(Y,\mu_Y)$ be $\sigma$-finite measure spaces. Let $T$ be a sublinear operator defined on simple functions supported on sets of finite measure in $X$, with values in measurable functions on $Y$, and suppose that it is of weak type $(p_0,p_0)$ and $(p_1,p_1)$, where $1\leq p_0<p_1\leq\infty$. If $p_0<p<p_1$, then $T$ extends to a bounded operator from $L^p(X)$ to $L^p(Y)$, and
\[
\|Tf\|_{L^p(Y)}\leq C\|f\|_{L^p(X)}.
\]
The constant depends only on $p_0,p,p_1$ and the weak-type constants of $T$.
In this result, sublinearity is understood with the following
bounds, almost everywhere on $Y$, for simple functions $f,g$
of finite measure support and scalars $\alpha$:
\[
 \begin{gathered}
 T0=0,\qquad |T(\alpha f)|=|\alpha|\,|Tf|,\\
 |T(f+g)|\leq|Tf|+|Tg|,
 \qquad |Tf-Tg|\leq|T(f-g)|.
 \end{gathered}
\]
The last bound holds for linear operators and for nonnegative
absolutely homogeneous sublinear operators; in the latter case,
apply subadditivity in both orders. The asserted extension
is unique among continuous extensions.
\end{theorem}

\begin{proof}
Start with a simple function $f$ of finite measure support.
First suppose $p_1<\infty$. Fix $\gamma>0$ and, for each $\lambda>0$, write
\[
f=f_0^\lambda+f_1^\lambda,
\qquad
f_0^\lambda:=f\mathbf1_{\{|f|>\gamma\lambda\}},
\qquad
f_1^\lambda:=f\mathbf1_{\{|f|\leq\gamma\lambda\}}.
\]
Sublinearity implies
\[
\{|Tf|>\lambda\}
\subseteq
\left\{|Tf_0^\lambda|>\frac{\lambda}{2}\right\}
\cup
\left\{|Tf_1^\lambda|>\frac{\lambda}{2}\right\}.
\]
If $C_0,C_1$ are the weak-type constants, we obtain
\[
\mu_Y(\{|Tf|>\lambda\})
\leq
C_0'\lambda^{-p_0}\|f_0^\lambda\|_{L^{p_0}(X)}^{p_0}
+
C_1'\lambda^{-p_1}\|f_1^\lambda\|_{L^{p_1}(X)}^{p_1}.
\]
Apply Cavalieri and Tonelli to the first term:
\[
\begin{aligned}
\int_0^\infty
\lambda^{p-p_0-1}
\int_{\{|f|>\gamma\lambda\}}|f(x)|^{p_0}\,d\mu_X(x)d\lambda
&=
\int_X|f(x)|^{p_0}
\int_0^{\frac{|f(x)|}{\gamma}}\lambda^{p-p_0-1}\,d\lambda d\mu_X(x)\\
&=
\frac{\gamma^{-(p-p_0)}}{p-p_0}\|f\|_{L^p(X)}^p.
\end{aligned}
\]
For the second term, since $p<p_1$,
\[
\begin{aligned}
\int_0^\infty
\lambda^{p-p_1-1}
\int_{\{|f|\leq\gamma\lambda\}}|f(x)|^{p_1}\,d\mu_X(x)d\lambda
&=
\int_X|f(x)|^{p_1}
\int_{\frac{|f(x)|}{\gamma}}^{\infty}\lambda^{p-p_1-1}\,d\lambda d\mu_X(x)\\
&=
\frac{\gamma^{p_1-p}}{p_1-p}\|f\|_{L^p(X)}^p.
\end{aligned}
\]
By Cavalieri, these two estimates imply $\|Tf\|_{L^p(Y)}^p\leq C\|f\|_{L^p(X)}^p$.

If $p_1=\infty$, the hypothesis at that endpoint means $\|Tg\|_{L^\infty(Y)}\leq C_1\|g\|_{L^\infty(X)}$. We may choose $C_1>0$ by increasing the bound if necessary.
Choose $\gamma=(2C_1)^{-1}$. Then
\[
\|Tf_1^\lambda\|_{L^\infty(Y)}
\leq
C_1\gamma\lambda
=\frac{\lambda}{2},
\]
so only the first term remains, which is integrated as above.
We have proved the strong bound for every simple function
of finite measure support.

Now let $f\in L^p(X)$. Since $p<\infty$, choose simple functions
of finite measure support $s_k$ with $s_k\to f$ in $L^p(X)$,
by Proposition~\ref{prop:densidad-funciones-simples-Lp}.
The difference bound and the estimate just obtained give
\[
 \|Ts_k-Ts_\ell\|_{L^p(Y)}
 \leq\|T(s_k-s_\ell)\|_{L^p(Y)}
 \leq C\|s_k-s_\ell\|_{L^p(X)}.
\]
By completeness, $Ts_k$ converges in $L^p(Y)$. If $u_k\to f$ is
another approximation, the same estimate applied to $s_k-u_k$ shows
that both sequences of images have the same limit.
Define $\widetilde Tf:=\lim_kTs_k$. Convergence of the norms gives
$\|\widetilde Tf\|_{L^p(Y)}\leq C\|f\|_{L^p(X)}$.
For $f,g\in L^p(X)$, taking approximations of both and passing
to the limit in the difference estimate gives
\[
 \|\widetilde Tf-\widetilde Tg\|_{L^p(Y)}
 \leq C\|f-g\|_{L^p(X)}.
\]
Thus the extension is continuous. It agrees with $T$ on simple
functions and is unique by density. If $T$ is linear, passing
to the limit in its linearity identities shows that $\widetilde T$ is linear as well.
\end{proof}

\begin{lemma}[Vitali selection for balls]\label{lem:seleccion-vitali-bolas}
Let $\mathcal B$ be a family of open balls in $\mathbb R^n$ with radii bounded above. There exists a countable disjoint subfamily $(B_j)_{j\in\mathbb N}\subseteq\mathcal B$ such that
\[
\bigcup_{B\in\mathcal B}B
\subseteq
\bigcup_{j\in\mathbb N}5B_j,
\]
where $5B_j$ denotes the concentric ball with five times the radius.
\end{lemma}

\begin{proof}
Let $R>0$ be an upper bound for the radii. For each $k\in\mathbb N_0$, consider the subfamily
\[
\mathcal B_k
:=
\left\{B\in\mathcal B\middle|2^{-k-1}R<r(B)\leq2^{-k}R\right\}.
\]
Proceed by levels. In $\mathcal B_0$, choose a maximal subfamily of disjoint balls. Once the families at levels $0,\dots,k-1$ have been chosen, choose in $\mathcal B_k$ a maximal subfamily of mutually disjoint balls that are also disjoint from all previously selected balls. Existence of each maximal family follows from Zorn's lemma. Since every family of disjoint open balls in $\mathbb R^n$ is countable---each ball contains a distinct point of $\mathbb Q^n$---the union of all selected families is countable; denote it by $(B_j)_{j\in\mathbb N}$.

Let $B\in\mathcal B_k$ be a ball that was not selected. By maximality, $B$ intersects a selected ball $B_j$ from the same or an earlier level. In both cases,
\[
r(B_j)\geq\frac{1}{2}r(B).
\]
If $x\in B$ and $c,c_j$ are the centers of $B$ and $B_j$, respectively, choose $z\in B\cap B_j$ and obtain
\[
\|x-c_j\|
\leq
\|x-c\|+\|c-z\|+\|z-c_j\|
<2r(B)+r(B_j)
\leq5r(B_j).
\]
Thus $B\subseteq5B_j$, proving that the dilates
$(5B_j)_{j\in\mathbb N}$ cover the union of the balls in $\mathcal B$.
\end{proof}

\begin{theorem}[Hardy--Littlewood maximal operator]\label{teo:maximal-hardy-littlewood}
For $f\in L^1_{\mathrm{loc}}(\mathbb R^n)$, define
\[
\mathcal Mf(x)
:=
\sup_{B\ni x}\frac{1}{|B|}\int_B|f(y)|\,dy,
\]
where the supremum is taken over all open balls $B\subseteq\mathbb R^n$ containing $x$. Then $\mathcal M$ is of weak type $(1,1)$: for $f\in L^1(\mathbb R^n)$ and $\lambda>0$,
\[
|\{x\in\mathbb R^n\mid \mathcal Mf(x)>\lambda\}|
\leq
\frac{C_n}{\lambda}\|f\|_{L^1(\mathbb R^n)}.
\]
Moreover, for $1<p\leq\infty$,
\[
\|\mathcal Mf\|_{L^p(\mathbb R^n)}
\leq
C_{n,p}\|f\|_{L^p(\mathbb R^n)}.
\]
\end{theorem}

\begin{proof}
For $f\in L^1_{\mathrm{loc}}$ and $\lambda>0$,
\[
 \{\mathcal Mf>\lambda\}
 =\bigcup_{\substack{B\text{ open ball}\\
                    |B|^{-1}\int_B|f|>\lambda}}B.
\]
This set is open; in particular, $\mathcal Mf$ is measurable,
with values allowed in $[0,\infty]$.
If $f\in L^1$, choose for each $x$ in this set a ball $B_x\ni x$
whose average exceeds $\lambda$. Then
$|B_x|<\displaystyle\frac{\|f\|_{L^1}}{\lambda}$, so their radii are uniformly bounded.
Lemma~\ref{lem:seleccion-vitali-bolas} gives disjoint balls $(B_j)$
such that $\{\mathcal Mf>\lambda\}\subseteq\bigcup_j5B_j$.
Therefore
\[
 |\{\mathcal Mf>\lambda\}|
 \leq5^n\sum_j|B_j|
 \leq\frac{5^n}{\lambda}\sum_j\int_{B_j}|f|
 \leq\frac{5^n}{\lambda}\|f\|_{L^1}.
\]
For $f\in L^\infty$, all averages are bounded by
$\|f\|_{L^\infty}$, proving the infinite endpoint.

On simple functions of finite measure support, $\mathcal M$
is nonnegative and absolutely homogeneous. The triangle
inequality for averages gives
$\mathcal M(f+g)\leq\mathcal Mf+\mathcal Mg$; applying it in both
orders also gives $|\mathcal Mf-\mathcal Mg|\leq\mathcal M(f-g)$.
Thus Theorem~\ref{teo:interpolacion-marcinkiewicz} yields, for
$1<p<\infty$,
\[
 \|\mathcal Ms\|_{L^p}\leq C_{n,p}\|s\|_{L^p}
\]
for every such simple function.

Now let $f\in L^p$. Hölder gives $f\in L^1_{\mathrm{loc}}$.
Choose simple functions $s_k$ of finite measure support with
$s_k\to f$ in $L^p$, as in
Proposition~\ref{prop:densidad-funciones-simples-Lp}.
For every ball $B$,
\[
 \frac{1}{|B|}\int_B\bigl||s_k|-|f|\bigr|
 \leq |B|^{-\frac{1}{p}}\|s_k-f\|_{L^p}\longrightarrow0.
\]
If $x\in B$, the average of $|f|$ on $B$ is therefore at most
$\liminf_k\mathcal Ms_k(x)$. Taking the supremum over balls
containing $x$ gives
$\mathcal Mf(x)\leq\liminf_k\mathcal Ms_k(x)$.
Fatou and convergence of the norms yield
\[
 \|\mathcal Mf\|_{L^p}^p
 \leq\liminf_k\|\mathcal Ms_k\|_{L^p}^p
 \leq C_{n,p}^p\|f\|_{L^p}^p.
\]
This proves the bound for the maximal operator defined by
averages and its finiteness almost everywhere.
\end{proof}

\begin{theorem}[Lebesgue differentiation]\label{teo:diferenciacion-lebesgue}
If $f\in L^1_{\mathrm{loc}}(\mathbb R^n)$, then, for almost every $x\in\mathbb R^n$,
\[
\lim_{r\to0^+}
\frac{1}{|B(x,r)|}
\int_{B(x,r)}|f(y)-f(x)|\,dy
=0.
\]
\end{theorem}

\begin{proof}
Fix $R>0$ and set $f_R:=f\mathbf1_{B(0,R+1)}\in L^1(\mathbb R^n)$. For $x\in B(0,R)$ and $0<r<1$, the averages of $f$ and $f_R$ over $B(x,r)$ agree. Define
\[
D h(x)
:=
\limsup_{r\to0^+}
\frac{1}{|B(x,r)|}\int_{B(x,r)}|h(y)-h(x)|\,dy.
\]
If $g\in C_c^\infty(\mathbb R^n)$, uniform continuity of $g$ gives $Dg(x)=0$ for every $x$. Moreover, for every $h\in L^1(\mathbb R^n)$,
\[
Dh(x)
\leq
\mathcal Mh(x)+|h(x)|
\]
for almost every $x$. Now choose $g\in C_c^\infty(\mathbb R^n)$ and use
\[
Df_R(x)
\leq
D(f_R-g)(x)+Dg(x)
\leq
\mathcal M(f_R-g)(x)+|f_R(x)-g(x)|.
\]
Thus, for $\varepsilon>0$,
\[
\{x\in B(0,R)\mid Df(x)>2\varepsilon\}
\subseteq
\{\mathcal M(f_R-g)>\varepsilon\}
\cup
\{|f_R-g|>\varepsilon\}.
\]
Weak type $(1,1)$ from Theorem~\ref{teo:maximal-hardy-littlewood} and Chebyshev give
\[
|\{x\in B(0,R)\mid Df(x)>2\varepsilon\}|
\leq
\frac{C_n}{\varepsilon}\|f_R-g\|_{L^1(\mathbb R^n)}.
\]
By Proposition~\ref{densidad de Cc en Lp}, $C_c^\infty(\mathbb R^n)$ is dense in $L^1(\mathbb R^n)$, so the right-hand side can be made arbitrarily small. Hence $Df=0$ almost everywhere on $B(0,R)$. Finally, let $R\to\infty$ through the integers.
\end{proof}

For Triebel--Lizorkin norms, controlling each block separately does not suffice: the $L^p(\mathbb R^n)$ norm is taken after the scales have been combined pointwise. We therefore need a vector-valued version of the maximal inequality. The following weighted result provides an efficient proof for part of the range of exponents and will also be useful in its own right.

\begin{lemma}[Weighted Fefferman--Stein inequality]
\label{lem:fefferman-stein-ponderada}
If $1<a<\infty$, $f\in L^1_{\mathrm{loc}}(\mathbb R^n)$, and $w\colon \mathbb R^n\longrightarrow[0,\infty]$ is measurable, then
\[
\int_{\mathbb R^n}(\mathcal Mf)^a w
\leq
C_{n,a}\int_{\mathbb R^n}|f|^a\mathcal Mw.
\]
The assertion is understood first for bounded compactly supported functions; the general case follows by truncation and monotone convergence.
\end{lemma}

\begin{proof}
It suffices to prove the estimate for a truncated version of the maximal operator considering only balls with radii in a compact interval of $(0,\infty)$; at the end, the truncations are removed by monotone convergence. Fix $\lambda>0$. Applying Lemma~\ref{lem:seleccion-vitali-bolas} to balls on which the average of $|f|$ exceeds $\lambda$ gives disjoint balls $(B_j)_{j\in\mathbb N}$ such that
\[
\{\mathcal Mf>\lambda\}
\subseteq
\bigcup_{j\in\mathbb N}5B_j.
\]
Since the contribution of $|f|\mathbf1_{\{|f|\leq\displaystyle\frac{\lambda}{2}\}}$ to the average over any ball is at most $\displaystyle\frac{\lambda}{2}$, for each $j\in\mathbb N$ we have
\[
|B_j|
\leq
\frac{2}{\lambda}
\int_{B_j\cap\left\{|f|>\frac{\lambda}{2}\right\}}|f(y)|\,dy.
\]
If $y\in B_j$, then $5B_j$ contains $y$, and the definition of $\mathcal Mw(y)$ gives
\[
\int_{5B_j}w
\leq
|5B_j|\,\mathcal Mw(y)
=5^n|B_j|\,\mathcal Mw(y).
\]
Taking the infimum over $y\in B_j$ and using the preceding inequality for $|B_j|$ yields
\[
\begin{aligned}
\int_{\{\mathcal Mf>\lambda\}}w
&\leq
\sum_{j\in\mathbb N}\int_{5B_j}w\\
&\leq
\frac{C_n}{\lambda}
\sum_{j\in\mathbb N}
\int_{B_j\cap\left\{|f|>\frac{\lambda}{2}\right\}}|f(y)|\mathcal Mw(y)\,dy\\
&\leq
\frac{C_n}{\lambda}
\int_{\left\{|f|>\frac{\lambda}{2}\right\}}|f|\mathcal Mw.
\end{aligned}
\]
Multiply by $a\lambda^{a-1}$ and integrate in $\lambda\in(0,\infty)$. By Tonelli,
\[
\begin{aligned}
\int_{\mathbb R^n}(\mathcal Mf)^aw
&=
a\int_0^\infty\lambda^{a-1}
\int_{\{\mathcal Mf>\lambda\}}w\,d\lambda\\
&\leq
C_{n,a}\int_{\mathbb R^n}|f(y)|\mathcal Mw(y)
\left(\int_0^{2|f(y)|}\lambda^{a-2}\,d\lambda\right)dy\\
&\leq
C_{n,a}\int_{\mathbb R^n}|f|^a\mathcal Mw.
\end{aligned}
\]
The inner integral is finite because $a>1$. The truncations are removed by Theorem~\ref{convergencia monotona}.
\end{proof}

\begin{theorem}[Vector-valued Fefferman--Stein inequality]\label{teo:fefferman-stein-vectorial}
Let $1<p,r<\infty$. If $J$ is a finite set and $(f_j)_{j\in J}$ is a family of measurable functions, then
\[
\left\|
\left(\sum_{j\in J}|\mathcal Mf_j|^r\right)^{\frac{1}{r}}
\right\|_{L^p(\mathbb R^n)}
\leq
C_{n,p,r}
\left\|
\left(\sum_{j\in J}|f_j|^r\right)^{\frac{1}{r}}
\right\|_{L^p(\mathbb R^n)}.
\]
The same inequality holds for countable families whenever the right-hand side is finite.
\end{theorem}

See also \cite[Section~2.2.2]{Triebel2}.

\begin{proof}
First suppose $p>r$. Then $\displaystyle\frac{p}{r}>1$. By duality of $L^{\displaystyle\frac{p}{r}}(\mathbb R^n)$,
\[
\left\|\sum_{j\in J}(\mathcal Mf_j)^r\right\|_{L^{\frac{p}{r}}(\mathbb R^n)}
=
\sup_{\substack{h\in L^{(\frac{p}{r})'}(\mathbb R^n),\ h\geq0\\\|h\|_{L^{(\frac{p}{r})'}(\mathbb R^n)}\leq1}}
\int_{\mathbb R^n}\sum_{j\in J}(\mathcal Mf_j)^r h.
\]
Apply Lemma~\ref{lem:fefferman-stein-ponderada} to each $f_j$, sum, and use Tonelli:
\[
\int_{\mathbb R^n}\sum_{j\in J}(\mathcal Mf_j)^r h
\leq
C\int_{\mathbb R^n}\sum_{j\in J}|f_j|^r\mathcal Mh.
\]
Hölder and continuity of $\mathcal M$ on $L^{(\displaystyle\frac{p}{r})'}(\mathbb R^n)$ give
\[
\int_{\mathbb R^n}\sum_{j\in J}|f_j|^r\mathcal Mh
\leq
C\left\|\sum_{j\in J}|f_j|^r\right\|_{L^{\frac{p}{r}}(\mathbb R^n)}.
\]
Taking the $r$th root gives the inequality when $p>r$.

Now suppose $p\leq r$ and fix a finite family $J$. The supremum defining $\mathcal M$ may be computed, up to a uniform dimensional constant, using a countable family of balls with rational centers and radii. For each $j\in J$ and the function $|f_j|$, measurably choose one of these balls $B_j(x)\ni x$ whose average is at least half the corresponding supremum. Once this choice has been made, define the linear operator
\[
T_jg(x):=\frac{1}{|B_j(x)|}\int_{B_j(x)}g(y)\,dy.
\]
The selection depends on the original family $(f_j)_{j\in J}$, but remains fixed while interpolating the diagonal operator
\[
\mathcal T(g_j)_{j\in J}:=(T_jg_j)_{j\in J}.
\]
For every $g_j$, we have $|T_jg_j|\leq\mathcal M g_j$. The scalar bound in Theorem~\ref{teo:maximal-hardy-littlewood} implies
\[
\begin{aligned}
\|\mathcal T(g_j)_{j\in J}\|_{L^p(\mathbb R^n;\ell^p(J))}^p
&=
\sum_{j\in J}\|T_jg_j\|_{L^p(\mathbb R^n)}^p\\
&\leq
C\sum_{j\in J}\|g_j\|_{L^p(\mathbb R^n)}^p.
\end{aligned}
\]
Moreover,
\[
\sup_{j\in J}|T_jg_j(x)|
\leq
\mathcal M\left(\sup_{j\in J}|g_j|\right)(x),
\]
so
\[
\|\mathcal T(g_j)_{j\in J}\|_{L^p(\mathbb R^n;\ell^\infty(J))}
\leq
C\|(g_j)_{j\in J}\|_{L^p(\mathbb R^n;\ell^\infty(J))}.
\]
If $r=p$, the desired estimate is precisely the first of these two bounds. Thus suppose $r>p$. There exists $0<\vartheta<1$ such that
\[
\frac{1}{r}=\frac{1-\vartheta}{p}.
\]
Interpolate the estimates with values in $\ell^p(J)$ and $\ell^\infty(J)$ using Theorem~\ref{teo:interpolacion-compleja-operadores} and Proposition~\ref{prop:interpolacion-ellp-valores-banach}. We obtain
\[
\|\mathcal T(g_j)_{j\in J}\|_{L^p(\mathbb R^n;\ell^r(J))}
\leq
C_{n,p,r}\|(g_j)_{j\in J}\|_{L^p(\mathbb R^n;\ell^r(J))}.
\]
Apply this bound to $g_j=|f_j|$. By the choice of $B_j(x)$,
\[
\mathcal Mf_j(x)
\leq
2T_j|f_j|(x),
\]
which proves the case $p\leq r$. The constant is independent of $J$. For a countable family, apply the result to the first $N$ coordinates and let $N\to\infty$ using monotone convergence.
\end{proof}

\begin{theorem}[Calderón--Zygmund decomposition]\label{teo:descomposicion-calderon-zygmund}
Let $f\in L^1(\mathbb R^n)$ and $\lambda>0$. There exist a countable family of disjoint half-open dyadic cubes $(Q_j)_{j\in\mathbb N}$ and functions $g,(b_j)_{j\in\mathbb N}$ such that
\[
f=g+\sum_{j\in\mathbb N}b_j
\quad\text{almost everywhere},
\]
\[
\|g\|_{L^\infty(\mathbb R^n)}\leq2^n\lambda,
\qquad
\|g\|_{L^1(\mathbb R^n)}\leq\|f\|_{L^1(\mathbb R^n)},
\]
\[
\operatorname{supp}(b_j)\subseteq \overline{Q_j},
\qquad
\int_{Q_j}b_j=0
\quad (j\in\mathbb N),
\]
and
\[
\sum_{j\in\mathbb N}|Q_j|
\leq
\lambda^{-1}\|f\|_{L^1(\mathbb R^n)},
\qquad
\sum_{j\in\mathbb N}\|b_j\|_{L^1(\mathbb R^n)}
\leq
2\|f\|_{L^1(\mathbb R^n)}.
\]
\end{theorem}

\begin{proof}
Use half-open dyadic cubes
$2^{-k}(m+[0,1)^n)$, with $k\in\mathbb Z$ and $m\in\mathbb Z^n$.
Any two are disjoint or one is contained in the other.
Consider the cubes whose average of $|f|$ exceeds $\lambda$.
If $x$ is a Lebesgue point of $f$ and $Q_k(x)$ is its cube of
side length $2^{-k}$, then
\[
 \frac{1}{|Q_k(x)|}\int_{Q_k(x)}|f(y)-f(x)|\,dy
 \leq C_n\frac{1}{|B(x,\sqrt n\,2^{-k})|}
                \int_{B(x,\sqrt n\,2^{-k})}|f(y)-f(x)|\,dy
 \longrightarrow0.
\]
The comparison holds because $Q_k(x)$ lies in that ball and
the ratio of their volumes is constant.
By Theorem~\ref{teo:diferenciacion-lebesgue}, almost every point
with $|f(x)|>\lambda$ belongs to some cube in the family.

Every cube in the family is contained in a maximal one.
Indeed, along its successive parents the side lengths tend
to infinity, and their averages are bounded by
$\displaystyle\frac{\|f\|_{L^1}}{|Q|}$, which tends to zero. The chain of cubes whose
averages exceed $\lambda$ must therefore end.
The maximal cubes $(Q_j)$ are countable and disjoint, and cover
$\{|f|>\lambda\}$ up to a null set. For the parent
$\widetilde Q_j$, maximality gives
\[
 \lambda<\frac{1}{|Q_j|}\int_{Q_j}|f|
 \leq2^n\frac{1}{|\widetilde Q_j|}\int_{\widetilde Q_j}|f|
 \leq2^n\lambda.
\]
Put $f_{Q_j}:=|Q_j|^{-1}\int_{Q_j}f$ and define
\[
 g:=f\mathbf1_{\mathbb R^n\setminus\bigcup_jQ_j}
          +\sum_j f_{Q_j}\mathbf1_{Q_j},
 \qquad
 b_j:=(f-f_{Q_j})\mathbf1_{Q_j}.
\]
Disjointness makes the sums pointwise unambiguous.
Outside the union, $|f|\leq\lambda$ almost everywhere, and on
$Q_j$, $|g|=|f_{Q_j}|\leq2^n\lambda$. Moreover,
\[
 \|g\|_{L^1}
 \leq\int_{\mathbb R^n\setminus\bigcup_jQ_j}|f|
          +\sum_j|f_{Q_j}|\,|Q_j|
 \leq\|f\|_{L^1}.
\]
Each $b_j$ vanishes outside $Q_j$, so its support lies in
$\overline{Q_j}$, and
$\int_{Q_j}b_j=\int_{Q_j}f-f_{Q_j}|Q_j|=0$.
Finally,
\[
 \sum_j|Q_j|\leq\lambda^{-1}\|f\|_{L^1},
 \qquad
 \sum_j\|b_j\|_{L^1}
 \leq\sum_j\left(\int_{Q_j}|f|+|f_{Q_j}|\,|Q_j|\right)
 \leq2\|f\|_{L^1}.
\]
The last bound also proves convergence of $\sum_jb_j$ in $L^1$;
disjointness gives $f=g+\sum_jb_j$ almost everywhere.
\end{proof}

The last result of this section directly connects frequency localization with the maximal operator. For a function whose spectrum is contained in a ball of radius comparable to $2^j$, the maximum of its translates, penalized at scale $2^{-j}$, is controlled by a Hardy--Littlewood maximal function of a power of the function. This is the form of Peetre's maximal inequality used in proofs of independence of the dyadic resolution and in almost diagonal arguments.

\begin{lemma}[Peetre maximal control]\label{lem:control-maximal-peetre}
Let $j\in\mathbb Z$, $0<a<\infty$, $N>\displaystyle\frac{n}{a}$, and $A>0$.
Suppose that $f\in\mathcal S'(\mathbb R^n)$, and identify it with its
canonical smooth representative. Suppose also that
\[
\operatorname{supp}(\widehat f)
\subseteq
\{\xi\in\mathbb R^n\mid \|\xi\|\leq A2^j\}.
\]
Then
\begin{equation}
\sup_{y\in\mathbb R^n}
\frac{|f(x-y)|}{(1+2^j\|y\|)^N}
\leq
C\bigl[\mathcal M(|f|^a)(x)\bigr]^{\frac{1}{a}},
\label{eq:peetre-maximal-control}
\end{equation}
where $C$ depends only on $n,a,N$ and $A$, at every point $x$ such that
$\mathcal M(|f|^a)(x)<\infty$. Consequently, if
\[
K_j(x):=2^{jn}K(2^jx)
\]
and
\[
|K(x)|\leq C_K(1+\|x\|)^{-N-n-1},
\]
then
\[
|K_j*f(x)|
\leq
C\bigl[\mathcal M(|f|^a)(x)\bigr]^{\frac{1}{a}}
\]
where $C$ also depends on $C_K$; at those points, the convolution converges
absolutely.
\end{lemma}

The maximal method for band-limited functions is developed in \cite[Sections~2.2--2.3]{Triebel2}.

\begin{proof}
Choose $\vartheta\in C_c^\infty(\mathbb R^n)$ such that $\vartheta=1$ on a neighborhood of the ball $\{\|\xi\|\leq A\}$, and set $\eta:=\mathcal F^{-1}\vartheta\in\mathcal S(\mathbb R^n)$. If $\eta_j(x):=2^{jn}\eta(2^jx)$, the convolution formula corresponding to the unitary Fourier normalization gives
\[
f=(2\pi)^{-\frac{n}{2}}\eta_j*f
\]
in $\mathcal S'(\mathbb R^n)$. The constant factor will be absorbed in the
estimates. In particular, by
Proposition~\ref{prop:regularizacion-temperada}, $f$ has the canonical
smooth representative of at most polynomial growth used in the
statement.

Since $\eta$ is a Schwartz function, for every $L>0$ there exists $C_L$ such that
\begin{equation}
|\eta_j(z)|
\leq
C_L2^{jn}(1+2^j\|z\|)^{-L}.
\label{eq:eta-schwartz-peetre}
\end{equation}

First suppose $a>1$. For $z\in\mathbb R^n$, the reproducing formula and the inequality
\[
1+2^j\|x-w\|
\leq
(1+2^j\|x-z\|)(1+2^j\|z-w\|)
\]
give, after dividing by $(1+2^j\|x-z\|)^N$,
\[
\frac{|f(z)|}{(1+2^j\|x-z\|)^N}
\leq
C_L\int_{\mathbb R^n}
2^{jn}(1+2^j\|z-w\|)^{-L+N}
\frac{|f(w)|}{(1+2^j\|x-w\|)^N}\,dw.
\]
Apply Hölder with exponents $a$ and $a'=\displaystyle\frac{a}{a-1}$. If $L$ is chosen so that $(L-N)a'>n$, the second factor remains uniformly bounded and we obtain
\[
\frac{|f(z)|}{(1+2^j\|x-z\|)^N}
\leq
C
\left(
\int_{\mathbb R^n}
2^{jn}
\frac{|f(w)|^a}{(1+2^j\|x-w\|)^{Na}}\,dw
\right)^{\frac{1}{a}}.
\]
Since $Na>n$, the decomposition of $\mathbb R^n$ into the ball $B(x,2^{-j})$ and the annuli
\[
2^{\ell-j}<\|x-w\|\leq2^{\ell+1-j},
\qquad \ell\in\mathbb N_0,
\]
proves
\begin{equation}
\int_{\mathbb R^n}
2^{jn}\frac{|f(w)|^a}{(1+2^j\|x-w\|)^{Na}}\,dw
\leq
C\mathcal M(|f|^a)(x),
\label{eq:nucleo-radial-control-maximal}
\end{equation}
since the contribution of the $\ell$th annulus is bounded by $C2^{-\ell(Na-n)}\mathcal M(|f|^a)(x)$. Taking the supremum over $z=x-y$ proves \eqref{eq:peetre-maximal-control} for $a>1$.

Now consider $0<a\leq1$. To avoid assuming in advance that the Peetre maximum is finite, fix $\varepsilon>0$ and define
\[
P_\varepsilon f(x)
:=
\sup_{z\in\mathbb R^n}
\frac{|f(z)|}{(1+2^j\|x-z\|)^N(1+\varepsilon2^j\|x-z\|)^M},
\]
where $0<\varepsilon\leq1$ and $M$ is chosen greater than a polynomial growth order of $f$. The regularity and growth established above imply $P_\varepsilon f(x)<\infty$. Using
\[
|f(w)|^{1-a}
\leq
P_\varepsilon f(x)^{1-a}
(1+2^j\|x-w\|)^{N(1-a)}
(1+\varepsilon2^j\|x-w\|)^{M(1-a)}
\]
in the preceding reproducing formula and applying the weight inequalities
\[
1+2^j\|x-w\|
\leq
(1+2^j\|x-z\|)(1+2^j\|z-w\|)
\]
and their analogue with $\varepsilon2^j$ gives, upon choosing $L$ sufficiently large in \eqref{eq:eta-schwartz-peetre},
\[
P_\varepsilon f(x)
\leq
C P_\varepsilon f(x)^{1-a}
\int_{\mathbb R^n}
2^{jn}\frac{|f(w)|^a}{(1+2^j\|x-w\|)^{Na}}\,dw.
\]
The constant $C$ is independent of $\varepsilon$. If $P_\varepsilon f(x)>0$, divide by $P_\varepsilon f(x)^{1-a}$; if it is zero, the conclusion is immediate. Using \eqref{eq:nucleo-radial-control-maximal} yields
\[
P_\varepsilon f(x)^a
\leq
C\mathcal M(|f|^a)(x).
\]
As $\varepsilon\to0^{+}$, the denominators decrease and
$P_\varepsilon f(x)$ increases to the unregularized Peetre maximum.
The preceding bound proves that this maximum, denoted by $Pf(x)$,
is finite. In this first step, the constant may depend on $M$,
chosen according to the growth of $f$.

Now repeat the estimate with the finite maximum $Pf(x)$ and
without the auxiliary weight. The reproducing formula and
\[
|f(w)|^{1-a}
\leq Pf(x)^{1-a}(1+2^j\|x-w\|)^{N(1-a)}
\]
give, for $z\in\mathbb R^n$,
\[
\frac{|f(z)|}{(1+2^j\|x-z\|)^N}
\leq C\,Pf(x)^{1-a}
\int_{\mathbb R^n}
2^{jn}\frac{|f(w)|^a}{(1+2^j\|x-w\|)^{Na}}\,dw.
\]
Indeed, use
$(1+2^j\|x-w\|)^N\leq
(1+2^j\|x-z\|)^N(1+2^j\|z-w\|)^N$
and choose $L>N$ in \eqref{eq:eta-schwartz-peetre}.
The remaining factor $(1+2^j\|z-w\|)^{-L+N}$ is at most one.
This time $L$ depends only on $N$, and no constant depends
on the growth order of $f$.
Take the supremum over $z\in\mathbb R^n$, divide by
$Pf(x)^{1-a}$ when it is positive, and apply
\eqref{eq:nucleo-radial-control-maximal}.
This gives \eqref{eq:peetre-maximal-control} with the dependence
of constants asserted in the statement.

For the second assertion, denote by $P_{j,N}f(x)$ the left-hand side of \eqref{eq:peetre-maximal-control}. Then
\[
\begin{aligned}
|K_j*f(x)|
&\leq
\int_{\mathbb R^n}|K_j(y)|\,|f(x-y)|\,dy\\
&\leq
P_{j,N}f(x)
\int_{\mathbb R^n}|K_j(y)|(1+2^j\|y\|)^N\,dy.
\end{aligned}
\]
The change of variables $z=2^jy$ and the hypothesis on $K$ show that the last integral is bounded uniformly in $j$. Finally, apply \eqref{eq:peetre-maximal-control}.
\end{proof}

\begin{proposition}[Continuity of translations]\label{prop:continuidad-traslaciones-Lp}
If $1\leq p<\infty$ and $f\in L^p(\mathbb R^n)$, then
\[
\|\tau_hf-f\|_{L^p(\mathbb R^n)}\longrightarrow0
\qquad\text{as }\|h\|\to0,
\quad
\tau_hf(x):=f(x+h).
\]
\end{proposition}

\begin{proof}
If $f\in C_c(\mathbb R^n)$, convergence follows from uniform continuity and the fact that all translates with $\|h\|\leq1$ are supported in a fixed compact set. For $f\in L^p(\mathbb R^n)$, choose $g\in C_c(\mathbb R^n)$ such that $\|f-g\|_{L^p(\mathbb R^n)}<\varepsilon$. Since translations are isometries of $L^p(\mathbb R^n)$,
\[
\|\tau_hf-f\|_{L^p(\mathbb R^n)}
\leq
2\|f-g\|_{L^p(\mathbb R^n)}
+
\|\tau_hg-g\|_{L^p(\mathbb R^n)}.
\]
The last term tends to zero as $\|h\|\to0$, after which we let $\varepsilon\to0^{+}$.
\end{proof}

\section{Euclidean mollifiers}\label{ap:molificadores-euclidianos}

The local regularization used in Sobolev spaces and partial differential equations is based on a single Euclidean construction. We formulate it with the properties needed later; see also \cite[Appendix C, \S5]{Evans}.

\begin{theorem}[Existence and properties of Euclidean mollifiers]\label{teo:molificadores-euclidianos}\index{Euclidean mollifier}\index{approximation of the identity}
There exists a function $\rho\in C_c^\infty(\mathbb R^n)$ such that $\rho\geq0$, $\operatorname{supp}(\rho)\subseteq\overline{B}_{\mathrm{euc}}(0,1)$, and $\displaystyle\int_{\mathbb R^n}\rho\,d\lambda_n=1$. For $\varepsilon>0$, define
\[
\rho_\varepsilon(x)
:=
\varepsilon^{-n}\rho\left(\frac{x}{\varepsilon}\right).
\]
Then $\rho_\varepsilon\geq0$, $\operatorname{supp}(\rho_\varepsilon)\subseteq\overline{B}_{\mathrm{euc}}(0,\varepsilon)$, and $\displaystyle\int_{\mathbb R^n}\rho_\varepsilon\,d\lambda_n=1$. Moreover, the following properties hold.
\begin{enumerate}[label=(\alph*)]
\item If $\Omega\subseteq\mathbb R^n$ is open, $u\in L^1_{\mathrm{loc}}(\Omega)$, and $\Omega_\varepsilon:=\{x\in\Omega\mid\operatorname{dist}_{\mathrm{euc}}(x,\mathbb R^n\setminus\Omega)>\varepsilon\}$, then
\[
u_\varepsilon(x)
:=
(u*\rho_\varepsilon)(x)
=
\int_{\mathbb R^n}u(x-y)\rho_\varepsilon(y)\,dy,
\qquad x\in\Omega_\varepsilon,
\]
is smooth and $D^\alpha u_\varepsilon=u*D^\alpha\rho_\varepsilon$ for every multi-index $\alpha$.
\item If $1\leq p<\infty$ and $u\in L^p_{\mathrm{loc}}(\Omega)$, then $u_\varepsilon\to u$ in $L^p(K)$ for every compact set $K\subseteq\Omega$ as $\varepsilon\to0^+$.
\item If $u$ is continuous, then $u_\varepsilon\to u$ uniformly on each compact subset of $\Omega$. If $u\colon \mathbb R^n\longrightarrow\mathbb R$ is uniformly continuous, convergence is uniform on $\mathbb R^n$.
\item If $u\colon \mathbb R^n\longrightarrow\mathbb R$ is $L$-Lipschitz, then $u_\varepsilon$ is also $L$-Lipschitz,
\[
|Du_\varepsilon(x)|\leq L
\]
for every $x\in\mathbb R^n$, and
\[
\|u_\varepsilon-u\|_{L^\infty(\mathbb R^n)}
\leq
L\varepsilon m_1(\rho),
\qquad
m_1(\rho):=\int_{\mathbb R^n}\|z\|\rho(z)\,dz.
\]
\end{enumerate}
\end{theorem}

\begin{proof}
Define
\[
\vartheta(x)
:=
\begin{cases}
\exp\left(-\frac{1}{1-\|x\|^2}\right),& \|x\|<1,\\
0,& \|x\|\geq1.
\end{cases}
\]
The function $\vartheta$ is smooth: inside the ball this is immediate, and as one approaches $\|x\|=1$, each derivative is the product of a rational function of $1-\|x\|^2$ and $\exp\left(-\displaystyle\frac{1}{1-\|x\|^2}\right)$, which tends to zero faster than any power of $1-\|x\|^2$. Since $\vartheta$ is nonnegative and not identically zero, the constant $\displaystyle c:=\left(\int_{\mathbb R^n}\vartheta\,d\lambda_n\right)^{-1}$ is well defined. The function $\rho:=c\vartheta$ has the required properties. The change of variables $y=\varepsilon z$ proves the assertions about $\rho_\varepsilon$.

Now let $u\in L^1_{\mathrm{loc}}(\Omega)$. If $x\in\Omega_\varepsilon$, the integral defining $u_\varepsilon(x)$ uses only values of $u$ in $\Omega$. It may also be written as
\[
u_\varepsilon(x)
=
\int_\Omega u(z)\rho_\varepsilon(x-z)\,dz.
\]
On a compact neighborhood of any point of $\Omega_\varepsilon$, the integrand vanishes outside a fixed compact subset of $\Omega$. Since $u$ is integrable on that compact set and all derivatives of $\rho_\varepsilon$ are bounded, dominated convergence allows differentiation under the integral sign any number of times. Thus $u_\varepsilon\in C^\infty(\Omega_\varepsilon)$ and $D^\alpha u_\varepsilon=u*D^\alpha\rho_\varepsilon$.

We now prove convergence in $L^p_{\mathrm{loc}}(X)$. First observe that, for every $v\in L^p(\mathbb R^n)$,
\[
\|v(\,\cdot-h)-v\|_{L^p(\mathbb R^n)}
\longrightarrow0
\qquad\text{as}\qquad h\to0.
\]
Indeed, by Proposition~\ref{densidad de Cc en Lp}, choose $\phi\in C_c^\infty(\mathbb R^n)$ with $\|v-\phi\|_{L^p(\mathbb R^n)}<\delta$. Translation invariance of Lebesgue measure gives
\[
\|v(\,\cdot-h)-v\|_{L^p(\mathbb R^n)}
\leq
2\delta+\|\phi(\,\cdot-h)-\phi\|_{L^p(\mathbb R^n)}.
\]
The last term converges to zero by uniform continuity of $\phi$ and because, for $\|h\|\leq1$, all functions involved are supported in a fixed compact set. Since $\delta>0$ is arbitrary, continuity of translations follows.

Let $K\subseteq\Omega$ be compact. Since $K$ and $\mathbb R^n\setminus\Omega$ are disjoint and the former is compact, the number $d_K:=\operatorname{dist}_{\mathrm{euc}}(K,\mathbb R^n\setminus\Omega)$ is positive. Set $\displaystyle\delta_K:=\displaystyle\frac{d_K}{3}$. Then the closed set $F_K:=\{x\in\mathbb R^n\mid\operatorname{dist}_{\mathrm{euc}}(x,K)\leq2\delta_K\}$ is contained in $\Omega$. Proposition~\ref{funciones flan}, applied to $F_K\subseteq\Omega$, provides $\eta\in C_c^\infty(\Omega)$ such that $0\leq\eta\leq1$ and $\eta=1$ on a neighborhood of $F_K$. Extend $v:=\eta u$ by zero to $\mathbb R^n$. If $0<\varepsilon<\delta_K$, then $u_\varepsilon=v*\rho_\varepsilon$ on $K$. Minkowski's integral inequality and the change of variables $y=\varepsilon z$ give
\[
\|v*\rho_\varepsilon-v\|_{L^p(\mathbb R^n)}
\leq
\int_{\mathbb R^n}\rho(z)
\|v(\,\cdot-\varepsilon z)-v\|_{L^p(\mathbb R^n)}\,dz.
\]
For each $z$, the integrand converges to zero and is dominated by $2\rho(z)\|v\|_{L^p(\mathbb R^n)}$. The dominated convergence theorem implies that $u_\varepsilon\to u$ in $L^p(K)$.

If $u$ is continuous, for $x\in\Omega_\varepsilon$ we have
\[
|u_\varepsilon(x)-u(x)|
\leq
\int_{\mathbb R^n}\rho_\varepsilon(y)|u(x-y)-u(x)|\,dy
\leq
\sup_{y\in\overline B_{\mathrm{euc}}(0,\varepsilon)}|u(x-y)-u(x)|.
\]
Uniform continuity of $u$ on a compact neighborhood of $K$ proves uniform convergence on $K$. The same argument applies on all of $\mathbb R^n$ when $u$ is uniformly continuous.

Finally, if $u$ is $L$-Lipschitz, then
\[
|u_\varepsilon(x)-u_\varepsilon(y)|
\leq
\int_{\mathbb R^n}\rho_\varepsilon(z)|u(x-z)-u(y-z)|\,dz
\leq
L\|x-y\|.
\]
Since $u_\varepsilon$ is smooth, taking difference quotients in unit directions gives $\|Du_\varepsilon\|\leq L$. Moreover,
\[
|u_\varepsilon(x)-u(x)|
\leq
\int_{\mathbb R^n}\rho(z)|u(x-\varepsilon z)-u(x)|\,dz
\leq
L\varepsilon\int_{\mathbb R^n}\|z\|\rho(z)\,dz,
\]
uniformly with respect to $x$. This completes the proof.
\end{proof}

\section{Topological vector spaces}\label{ap:evt}

The function spaces arising in global analysis are often not naturally described by a single norm. This occurs, for example, with spaces of compactly supported smooth functions, spaces of test sections, and various weak topologies associated with duality problems. We therefore recall the language of topological vector spaces, which allows the linear and topological structures of a space to be studied together.

\begin{definition}\label{def: espacio vectorial topologico}\index{topological vector space}
Let $X$ be a vector space over a field $\mathbb K=\mathbb R,\mathbb C$ equipped with a topology $\tau$. The pair $(X,\tau)$ is called a \textit{topological vector space} if the functions $X\times X\longrightarrow X,\quad (x,y)\mapsto x+y$ and $\mathbb K\times X\longrightarrow X,\quad (\lambda,x)\mapsto \lambda x$ are continuous.
\end{definition}

Continuity of the linear operations has immediate consequences. For example, translations and nonzero dilations are homeomorphisms. In particular, the topology of a topological vector space is determined by its behavior around the origin.

\begin{proposition}\label{prop: topologia-determinada-por-base-local}
Let $X$ be a topological vector space, and let $\mathcal B$ be a neighborhood basis at the origin. Then, for each $x\in X$, the family
\[
x+\mathcal B:=\{x+U\mid U\in\mathcal B\}
\]
is a neighborhood basis at $x$. In particular, the family
\[
\{x+U\mid x\in X,\ U\in\mathcal B\}
\]
determines the open sets of $X$: a set $O\subseteq X$ is open if and only if, for each $x\in O$, there exists $U\in\mathcal B$ such that
$x+U\subseteq O$. If the elements of $\mathcal B$ are open, their translates form a basis of open sets. Consequently, the topology of a topological vector space is uniquely determined by any neighborhood basis at the origin.
\end{proposition}

\begin{proof}
Let $x\in X$. Since the translation $\tau_x\colon X\longrightarrow X$, given by $\tau_x(y)=x+y$, is a homeomorphism, $V$ is a neighborhood of $0$ if and only if $x+V$ is a neighborhood of $x$. Thus, if $\mathcal B$ is a local basis at $0$, then $x+\mathcal B$ is a local basis at $x$. If $O$ is open and $x\in O$, there exists $U\in\mathcal B$ such that $x+U\subseteq O$. Conversely, if this condition holds at each $x\in O$, every point of $O$ has a neighborhood contained in $O$ and is an interior point. This characterizes the open sets and proves uniqueness of the topology. When the elements of $\mathcal B$ are open, so are their translates, which therefore form a basis of open sets.
\end{proof}

To describe useful local bases, we need a few elementary geometric notions in vector spaces.

\begin{definition}[Absorbing set]\label{def: conjunto absorbente}\index{absorbing set}
Let $X$ be a vector space over $\mathbb R$ or $\mathbb C$. A set $A \subset X$ is said to be \emph{absorbing} if, for every $x \in X$, there exists $\lambda > 0$ such that $x \in tA$, $\forall t \geq \lambda$, where $tA := \{ta \mid a \in A\}$.
\end{definition}

\begin{definition}[Balanced set]\label{def: conjunto balanceado}\index{balanced set}
Let $X$ be a vector space over $\mathbb R$ or $\mathbb C$. A set $B \subset X$ is said to be \emph{balanced} if, for every $x \in B$ and every scalar $\alpha$ with $|\alpha|\leq 1$, we have $\alpha x \in B$.

\end{definition}

\begin{definition}[Absolutely convex set]\label{def: conjunto absolutamente convexo}\index{absolutely convex set}
Let $X$ be a vector space over $\mathbb R$ or $\mathbb C$. A subset $A\subseteq X$ is called \textit{absolutely convex} if it is convex and balanced.
\end{definition}

As a first example, every neighborhood of the origin absorbs the entire space.

\begin{lemma}\label{lem: vecindades-absorbentes}
Every neighborhood of the origin in a topological vector space is absorbing.
\end{lemma}

\begin{proof}
Let $U$ be a neighborhood of the origin, and let $x\in X$. Since the map
\[
\mathbb K\longrightarrow X,
\qquad
\lambda\longmapsto \lambda x,
\]
is continuous and $\lambda x\to 0$ as $\lambda\to 0$, there exists $\delta>0$ such that $\lambda x\in U$ if $|\lambda|<\delta$. Taking $t>\delta^{-1}$ gives $t^{-1}x\in U$ and hence $x\in tU$. Since $x\in X$ was arbitrary, $U$ is absorbing.
\end{proof}

\begin{lemma}\label{lem: vecindades pequenas en EVT}
Let $X$ be a topological vector space. Then:
\begin{enumerate}[label=(\alph*)]
 \item for every neighborhood $U$ of the origin, there is a neighborhood $V$ of the origin such that $V+V\subseteq U$;
 \item for every neighborhood $U$ of the origin, there is a balanced neighborhood $B$ of the origin such that $B\subseteq U$;
 \item if $X$ is Hausdorff, then, for every neighborhood $U$ of the origin, there is a neighborhood $W$ of the origin such that $\overline W\subseteq U$.
\end{enumerate}
\end{lemma}

\begin{proof}
For $(a)$, use continuity of addition at $(0,0)$. Since $U$ is a neighborhood of $0$, there are neighborhoods $V_1,V_2$ of the origin such that $V_1+V_2\subseteq U$. Taking $V=V_1\cap V_2$ gives $V+V\subseteq U$.

For $(b)$, use continuity of scalar multiplication at $(0,0)$. There exist $\varepsilon>0$ and a neighborhood $V$ of the origin such that $\lambda V\subseteq U$ if $|\lambda|<\varepsilon$. Taking $0<r<\varepsilon$ and defining
\[
B:=\bigcup_{|\lambda|\leq r}\lambda V,
\]
we find that $B$ is a neighborhood of the origin, $B\subseteq U$, and $B$ is balanced.

For $(c)$, by part (a), choose a neighborhood \(V_0\) of the origin such that \(V_0+V_0\subseteq U\). The set
\[
V:=V_0\cap(-V_0)
\]
is a symmetric neighborhood of the origin and satisfies \(V+V\subseteq U\). Let \(W:=\operatorname{Int}(V)\), which is an open symmetric neighborhood of the origin. If \(x\in\overline W\), the neighborhood \(x+W\) of \(x\) intersects \(W\). Thus there exist \(w_1,w_2\in W\) such that \(x+w_1=w_2\), and
\[
x=w_2-w_1\in W+W\subseteq V+V\subseteq U.
\]
Hence \(\overline W\subseteq U\). In particular, this argument directly proves regularity of every Hausdorff topological vector space.
\end{proof}

\begin{remark}\label{rem: base-local-origen}
Proposition~\ref{prop: topologia-determinada-por-base-local} shows that, in a topological vector space, it suffices to describe a local basis at the origin. Lemma~\ref{lem: vecindades pequenas en EVT} allows neighborhoods of the origin to be refined so that they have useful properties under addition and scalar multiplication; under the Hausdorff assumption, we can find a basis of closed neighborhoods.
\end{remark}

In a topological vector space, having a countable local basis at the origin is equivalent to being first countable. This property gives a characterization of metrizable topological vector spaces.

\begin{theorem}\label{teo: criterio de metrizabilidad para EVT}\index{metrizability criterion for TVSs}
A topological vector space $X$ is metrizable if and only if:
\begin{enumerate}[label=(\alph*)]
 \item $X$ is Hausdorff;
 \item $X$ has a countable local basis at $0$.
\end{enumerate}
\end{theorem}

The construction of a translation-invariant metric from a countable local basis, together with the immediate converse, can be found in \cite{Rudin1991,Treves1967}.

In metrizable spaces, continuity can be proved using sequences. This observation will be useful when working with Fréchet spaces and other metrizable spaces.

\begin{proposition}\label{prop: continuidad secuencial en espacios metrizables}
Let $X$ be a metrizable topological vector space, and let $Y$ be a topological vector space. A linear operator $L\colon X\longrightarrow Y$ is continuous if and only if it is sequentially continuous at the origin; that is, if, for every sequence $\{x_n\}_{n\in\mathbb N}$ in $X$ such that $x_n\to 0$, we have $L(x_n)\to 0$ in $Y$.
\end{proposition}

\begin{proof}
If $L$ is continuous and $x_n\to 0$ in $X$, continuity gives $L(x_n)\to L(0)=0$ in $Y$.

Conversely, suppose that $L$ is sequentially continuous at the origin. We prove that $L$ is continuous. Since $L$ is linear, it suffices to prove continuity at $0$. Suppose, for a contradiction, that $L$ is not continuous at $0$. There is then a neighborhood $U$ of the origin in $Y$ such that, for every neighborhood $V$ of the origin in $X$, we have $L(V)\not\subseteq U$.

Since $X$ is metrizable, it is first countable. Let $\{W_n\}_{n\in\mathbb N}$ be a countable local basis at the origin in $X$. Define
$
V_n:=\displaystyle\bigcap_{k=1}^{n}W_k,
$
for each $n\in\mathbb N$. Then each $V_n$ is a neighborhood of the origin, and the family $\{V_n\}_{n\in\mathbb N}$ is decreasing and remains a local basis at the origin.

By the choice of $U$, for each $n\in\mathbb N$ we can choose $x_n\in V_n$ such that $L(x_n)\notin U$.

Let us show that $x_n\to 0$ in $X$. Let $W$ be an arbitrary neighborhood of the origin. Since $\{V_n\}_{n\in\mathbb N}$ is a local basis at the origin, there exists $N\in\mathbb N$ such that $V_N\subseteq W$. Because the family is decreasing, $V_n\subseteq V_N\subseteq W$ for every $n\geq N$. Consequently, $x_n\in W$ for every $n\geq N$. This proves that $x_n\to 0$ in $X$.

By sequential continuity at the origin, $L(x_n)\to 0$ in $Y$. This is impossible, however, because $U$ is a neighborhood of the origin in $Y$ and $L(x_n)\notin U$ for every $n\in\mathbb N$. This contradiction shows that $L$ is continuous at $0$. Finally, linearity implies that $L$ is continuous.
\end{proof}

Boundedness in topological vector spaces depends on how sets are absorbed by neighborhoods of the origin, rather than on a metric.

\begin{definition}\label{def: conjunto topologicamente acotado}\index{topologically bounded}
Let $X$ be a topological vector space. A set $E\subseteq X$ is said to be \textit{topologically bounded} if, for each neighborhood $U$ of $0$, there exists $t>0$ such that $E\subseteq tU$.
\end{definition}

When the topology $\tau$ is generated by a metric $d$, metric boundedness and topological boundedness may differ. The metric $d_{1}=\displaystyle\frac{d}{d+1}$ generates the same topology as $d$, but, since $d_{1}\leq 1$, every set in $X$ is bounded with respect to $d_{1}$. Nevertheless, not every subset of $X$ need be topologically bounded.

We will also need the sequential notion of a Cauchy sequence. In nonmetrizable topological vector spaces, general completeness is formulated using Cauchy filters or nets. In this appendix, we use the sequential version only when it suffices, particularly in metrizable spaces.

\begin{definition}\label{def: sucesion de Cauchy EVT}\index{sequentially complete}
Let $X$ be a topological vector space. A sequence $\{x_{n}\}$ in $X$ is called a Cauchy sequence if, for every neighborhood $U$ of the origin, there exists $N\in \mathbb N$ such that
\[
x_{n}-x_{m}\in U,\quad \forall n,m\geq N.
\]
We say that $X$ is \textit{sequentially complete} if every Cauchy sequence in $X$ converges.
\end{definition}

\begin{proposition}\label{prop: sucesion de Cauchy tiene rango topologicamente acotado}
Let $X$ be a topological vector space, and let $\{x_{n}\}$ be a Cauchy sequence in $X$. Then the range of the sequence $\{x_{n}\mid n\in \mathbb N\}$ is topologically bounded.
\end{proposition}

\begin{proof}
Let $U$ be a neighborhood of the origin. By Lemma~\ref{lem: vecindades pequenas en EVT}, there is a balanced neighborhood $W$ of the origin such that $W+W\subseteq U$. Since $\{x_n\}$ is Cauchy, there exists $N\in\mathbb N$ such that $x_n-x_N\in W$ for every $n\geq N$. Moreover, since $W$ is absorbing, there exists $t_0>0$ such that $x_N\in tW$ for every $t\geq t_0$. If $t\geq \displaystyle\max\{1,t_0\}$ and $n\geq N$, then $x_n=x_N+(x_n-x_N)\in tW+W\subseteq t(W+W)\subseteq tU$.

On the other hand, the finite set $\{x_1,\dots,x_{N-1}\}$ is topologically bounded because each of its points is absorbed by $U$. Increasing $t$ if necessary, we find that all terms of the sequence belong to a single multiple of $U$. Since $U$ was arbitrary, the range of the sequence is topologically bounded.
\end{proof}

The most important class of topological vector spaces for our purposes will be that of locally convex spaces. Their importance stems from the fact that their topologies can be described by seminorms, which is particularly suitable for spaces of functions and sections.

\begin{definition}\label{def: espacio localmente convexo}\index{locally convex space}
A topological vector space $X$ is said to be locally convex if it has a local basis at $0$ consisting of convex sets; equivalently, every neighborhood of $0$ contains a convex neighborhood of $0$.
\end{definition}

\begin{proposition}\label{prop: localmente convexo tiene una base de vecindades absolutamente convexas}
A locally convex topological vector space admits a local basis at $0$ consisting of balanced convex sets.
\end{proposition}

\begin{proof}
Let $U$ be a neighborhood of the origin. Since $X$ is locally convex, there is a convex neighborhood $V$ of the origin such that $V\subseteq U$. By Lemma~\ref{lem: vecindades pequenas en EVT}, there is a balanced neighborhood $B$ of the origin such that $B\subseteq V$. The convex hull of $B$, denoted by $\operatorname{conv}(B)$, is then balanced and convex. Moreover, since $V$ is convex and $B\subseteq V$, we have $\operatorname{conv}(B)\subseteq V\subseteq U$. Thus every neighborhood of the origin contains a balanced convex neighborhood of the origin.
\end{proof}

\begin{definition}\label{def: seminorma}\index{seminorm}
Let $X$ be a vector space over $\mathbb K=\mathbb R$ or $\mathbb C$. A function $p\colon X\longrightarrow [0,\infty)$ is a seminorm if
\[
p(x+y)\leq p(x)+p(y)
\]
for any $x,y\in X$ and
\[
p(tx)=|t|p(x)
\]
for every $x\in X$ and any $t\in \mathbb K$.
\end{definition}

Every absorbing absolutely convex set has an associated seminorm given by its Minkowski functional.

\begin{proposition}\label{prop: minkowski-seminorma}
Let $X$ be a vector space over $\mathbb K$, where $\mathbb K=\mathbb R$ or $\mathbb C$, and let $A\subseteq X$ be an absorbing set. If $A$ is absolutely convex, its Minkowski functional
\[
p_A(x):=\inf\{\lambda>0 \mid x\in \lambda A\},
\qquad x\in X,
\]
is a seminorm on $X$.
\end{proposition}

\begin{proof}
Since $A$ is absorbing, for each $x\in X$ there exists $\lambda>0$ such that $x\in \lambda A$. Thus the set in the definition of $p_A(x)$ is nonempty, and consequently $p_A(x)$ is well defined for every $x\in X$.

First consider homogeneity. Let $\alpha\in\mathbb K$ and $x\in X$. If $\alpha=0$, then $p_A(\alpha x)=p_A(0)=0$. Now suppose that $\alpha\neq 0$. Write $\alpha=|\alpha|\theta$, where $|\theta|=1$. Since $A$ is balanced, $\theta A=A$. Therefore, for each $\lambda>0$,
\[
\alpha x\in \lambda A
\quad\Longleftrightarrow\quad
|\alpha|\theta x\in \lambda A
\quad\Longleftrightarrow\quad
x\in \frac{\lambda}{|\alpha|}A.
\]
It follows that
\[
p_A(\alpha x)=|\alpha|\,p_A(x).
\]

We now prove subadditivity. Let $x,y\in X$ and $\varepsilon>0$. By the definition of the infimum, there exist $\lambda,\mu>0$ such that
\[
x\in \lambda A,\qquad
y\in \mu A,\qquad
\lambda<p_A(x)+\varepsilon,\qquad
\mu<p_A(y)+\varepsilon.
\]
Since $A$ is convex, we have
\[
\lambda A+\mu A\subseteq (\lambda+\mu)A.
\]
If $a,b\in A$, then
\[
\lambda a+\mu b
=
(\lambda+\mu)
\left(
\frac{\lambda}{\lambda+\mu}a+
\frac{\mu}{\lambda+\mu}b
\right),
\]
and the term in parentheses belongs to $A$ by convexity. Therefore,
\[
x+y\in (\lambda+\mu)A,
\]
whence
\[
p_A(x+y)\leq \lambda+\mu
<
p_A(x)+p_A(y)+2\varepsilon.
\]
Since $\varepsilon>0$ was arbitrary, we conclude that
\[
p_A(x+y)\leq p_A(x)+p_A(y).
\]

Thus $p_A$ is a seminorm on $X$.
\end{proof}

\begin{definition}\label{def: topologia inducida por seminormas}\index{topology induced by seminorms}
Let $X$ be a vector space, and let $\mathcal P$ be a family of seminorms on $X$. The \textit{topology induced by $\mathcal P$} is the topology whose local basis at the origin consists of all finite intersections of sets of the form
\[
\{x\in X\mid p(x)<r\},
\qquad p\in\mathcal P,\quad r>0.
\]
Equivalently, a neighborhood basis at a point $x\in X$ is given by sets of the form
\[
x+\bigcap_{i=1}^{k}\{y\in X\mid p_i(y)<r_i\},
\]
where $p_1,\dots,p_k\in\mathcal P$ and $r_1,\dots,r_k>0$.
\end{definition}

\begin{theorem}\label{seminormas topologia}\index{topology induced by seminorms}\label{teo: seminormas topologia}
Let $X$ be a locally convex topological vector space, and let $\mathcal F$ be a local basis at the origin consisting of open absolutely convex neighborhoods. Then the family $\{p_U\mid U\in \mathcal F\}$ of Minkowski functionals is a family of continuous seminorms generating the topology of $X$.

Conversely, given a family $\mathcal P$ of seminorms on a vector space $X$, the topology induced by $\mathcal P$ makes $X$ a locally convex topological vector space, and every seminorm $p\in\mathcal P$ is continuous in this topology.
\end{theorem}

\begin{proof}
First suppose that $\mathcal F$ is a local basis at the origin consisting of open absolutely convex neighborhoods. Let $U\in\mathcal F$. Since $U$ is a neighborhood of the origin in a topological vector space, Lemma~\ref{lem: vecindades-absorbentes} shows that $U$ is absorbing. Proposition~\ref{prop: minkowski-seminorma} therefore ensures that its Minkowski functional $p_U$ is a seminorm on $X$.

We show that $p_U$ is continuous. Since $p_U$ is a seminorm, it suffices to prove continuity at the origin. Let $\varepsilon>0$. Multiplication by the nonzero scalar $\displaystyle\frac{\varepsilon}{2}$ is a homeomorphism of $X$ onto itself and $U$ is a neighborhood of the origin, so the set $\displaystyle\frac{\varepsilon}{2}U$ is also a neighborhood of the origin. If $x\in \displaystyle\frac{\varepsilon}{2}U$, then
\[
p_U(x)\leq \frac{\varepsilon}{2}<\varepsilon.
\]
Consequently,
\[
\frac{\varepsilon}{2}U\subseteq \{x\in X\mid p_U(x)<\varepsilon\},
\]
which proves that $p_U$ is continuous at $0$. Finally, for any $x,y\in X$, subadditivity of $p_U$ gives $p_U(x)\leq p_U(x-y)+p_U(y)$ and $p_U(y)\leq p_U(x-y)+p_U(x)$, whence $|p_U(x)-p_U(y)|\leq p_U(x-y)$.

Thus continuity of $p_U$ at the origin implies continuity at every point of $X$.

We now show that the family $\{p_U\mid U\in\mathcal F\}$ generates the original topology. For each $U\in\mathcal F$, since $U$ is open, balanced, and convex, we have
\[
U=\{x\in X\mid p_U(x)<1\}.
\]
If $p_U(x)<1$, there exists $\lambda<1$ such that $x\in \lambda U$, and balancedness of $U$ gives $x\in U$. Conversely, if $x\in U$, openness of $U$ and continuity of the map $t\mapsto tx$ give $r>1$ such that $rx\in U$. Therefore $x\in \displaystyle\frac{1}{r} U$ and $p_U(x)\leq \displaystyle\frac{1}{r}<1$. Thus the sets defined by the seminorms $p_U$ reproduce the local basis $\mathcal F$ and, by Proposition~\ref{prop: topologia-determinada-por-base-local}, generate the topology of $X$.

We now prove the converse assertion. Let $\mathcal P$ be a family of seminorms on the vector space $X$, and denote by $\mathcal B$ the collection of all finite intersections of sets of the form
\[
\{x\in X\mid p(x)<r\},
\qquad p\in\mathcal P,\quad r>0.
\]
Each of these sets is balanced and convex, as are their finite intersections. Define a topology $\tau$ by declaring $O\subseteq X$ open if, for each $x\in O$, there exists $B\in\mathcal B$ such that $x+B\subseteq O$. This makes $\mathcal B$ a local basis at the origin.

Addition is continuous at $(0,0)$ because, given a basic element
\[
B=\bigcap_{i=1}^{k}\{x\in X\mid p_i(x)<r_i\},
\]
the set
\[
C=\bigcap_{i=1}^{k}\left\{x\in X\mid p_i(x)<\frac{r_i}{2}\right\}
\]
satisfies $C+C\subseteq B$. Continuity of scalar multiplication follows similarly: if $\alpha_0\in\mathbb K$, $x_0\in X$, and $B$ is as above, choose $\delta>0$ and a basic neighborhood $C$ of $0$ such that
\[
|\alpha-\alpha_0|p_i(x_0)<\frac{r_i}{2}
\quad\text{and}\quad
|\alpha|p_i(h)<\frac{r_i}{2}
\]
for $|\alpha-\alpha_0|<\delta$, $h\in C$, and every $i\in\{1,\dots,k\}$. Then $\alpha(x_0+h)-\alpha_0x_0\in B$. Thus $\tau$ makes $X$ a topological vector space. Since $\mathcal B$ consists of convex sets, the space is locally convex. Finally, each seminorm $p\in\mathcal P$ is continuous because $\{x\in X\mid p(x)<\varepsilon\}$ is a neighborhood of the origin for every $\varepsilon>0$.
\end{proof}

We now give some consequences of the description by seminorms.

\begin{proposition}\label{prop: Hausdorff si singleton cerrado}
Let $X$ be a topological vector space. Then $X$ is Hausdorff if and only if the set $\{0\}$ is closed.
\end{proposition}

\begin{proof}
First suppose that $X$ is Hausdorff. Since every Hausdorff space is $T_1$, its singleton sets are closed. In particular, $\{0\}$ is closed.

Conversely, suppose that $\{0\}$ is closed. Let $x,y\in X$ with $x\neq y$. Then $x-y\neq0$ and hence $x-y\in X\setminus\{0\}$. Since $X\setminus\{0\}$ is open, there is a neighborhood $U$ of the origin such that
\[
x-y+U\subseteq X\setminus\{0\}.
\]
In particular, $x-y\notin -U$. Replacing $U$ by $U\cap(-U)$, we may assume that $U$ is symmetric, that is, $U=-U$. Thus $x-y\notin U$.

Since addition is continuous at $(0,0)$, there is a neighborhood $V$ of the origin such that $V+V\subseteq U$. Replacing $V$ by $V\cap(-V)$ again, we may assume that $V$ is symmetric.

Consider the neighborhoods $x+V$ and $y+V$. We claim that they are disjoint. If there were $z\in(x+V)\cap(y+V)$, there would be $v_1,v_2\in V$ such that $z=x+v_1=y+v_2$. Therefore $x-y=v_2-v_1$. Since $V$ is symmetric, $-v_1\in V$, and then $v_2-v_1\in V+V\subseteq U$, contradicting $x-y\notin U$. Therefore $(x+V)\cap(y+V)=\varnothing$.

We have found disjoint neighborhoods of $x$ and $y$. Since $x$ and $y$ were arbitrary, $X$ is Hausdorff.
\end{proof}

\begin{proposition}\label{prop: Hausdorff seminormas separan puntos}
Let $X$ be a vector space, and let $\mathcal P$ be a family of seminorms on $X$. The topology induced by $\mathcal P$ is Hausdorff if and only if $\mathcal P$ separates points, that is, if, for every $x\in X$ with $x\neq0$, there exists $p\in\mathcal P$ such that $p(x)\neq0$. Equivalently,
\[
\bigcap_{p\in\mathcal P}\{x\in X\mid p(x)=0\}=\{0\}.
\]
\end{proposition}

\begin{proof}
First suppose that the topology induced by $\mathcal P$ is Hausdorff. Let $x\in X$ with $x\neq0$. Since the topology is Hausdorff, there are disjoint neighborhoods of the origin and of $x$. In particular, there is a neighborhood $U$ of the origin such that $x\notin U$. By the definition of the induced topology, there exist $p_1,\dots,p_k\in\mathcal P$ and $r_1,\dots,r_k>0$ such that
\[
\bigcap_{i=1}^{k}\{y\in X\mid p_i(y)<r_i\}\subseteq U.
\]
Since $x\notin U$, we also have
\[
x\notin\bigcap_{i=1}^{k}\{y\in X\mid p_i(y)<r_i\}.
\]
Therefore there exists $i\in\{1,\dots,k\}$ such that $p_i(x)\geq r_i>0$. In particular, $p_i(x)\neq0$. Thus the family $\mathcal P$ separates points.

Conversely, suppose that $\mathcal P$ separates points. Consider
\[
N:=\bigcap_{p\in\mathcal P}\{x\in X\mid p(x)=0\}.
\]
Since every seminorm satisfies $p(0)=0$, we have $0\in N$. On the other hand, if $x\neq0$, the hypothesis gives $p\in\mathcal P$ such that $p(x)\neq0$, so $x\notin\{y\in X\mid p(y)=0\}$. Consequently, $x\notin N$. This proves that $N=\{0\}$.

We now show that $N$ is closed. For each $p\in\mathcal P$, the set $\{x\in X\mid p(x)=0\}$ is closed in the induced topology. If $x_0$ does not belong to this set, then $p(x_0)>0$. Taking $r:=\displaystyle\frac{p(x_0)}{2}$, the neighborhood
\[
x_0+\{y\in X\mid p(y)<r\}
\]
of $x_0$ does not intersect $\{x\in X\mid p(x)=0\}$. If $p(x)=0$ and $p(x-x_0)<r$, then
\[
p(x_0)=p(x_0-x+x)\leq p(x_0-x)+p(x)<r,
\]
contradicting the choice of $r$. Therefore $\{x\in X\mid p(x)=0\}$ is closed.

Since $N$ is an intersection of closed sets, $N$ is closed. But $N=\{0\}$, so $\{0\}$ is closed. By Proposition~\ref{prop: Hausdorff si singleton cerrado}, the topology induced by $\mathcal P$ is Hausdorff.
\end{proof}

\begin{corollary}\label{cor: topologia generada por seminormas hausdorff y topologicamente acotado}
Let $\mathcal P$ be a family of seminorms on a vector space $X$, and let $\tau_{\mathcal P}$ be the locally convex topology generated by $\mathcal P$. Then:
\begin{enumerate}[label=(\alph*)]
 \item $\tau_{\mathcal P}$ is Hausdorff if and only if $p(x)=0$ for every $p\in \mathcal P$ implies $x=0$.
 \item A set $A\subseteq X$ is topologically bounded if and only if the set $p(A)\subseteq \mathbb R$ is bounded in $\mathbb R$ for every $p\in \mathcal P$.
\end{enumerate}
A topological vector space is called \textit{seminormable} if its topology can be induced by a seminorm $p$. The space $(X,p)$ is called a \textit{seminormed} space.
\end{corollary}

\begin{proof}
The first assertion is Proposition~\ref{prop: Hausdorff seminormas separan puntos}.

For $(b)$, first suppose that $A$ is topologically bounded. If $p\in\mathcal P$, then $\{x\in X\mid p(x)<1\}$ is a neighborhood of the origin. There exists $t>0$ such that $A\subseteq t\{x\in X\mid p(x)<1\}$. Therefore, if $a\in A$, we have $p(a)<t$, so $p(A)$ is bounded.

Conversely, suppose that $p(A)$ is bounded for every $p\in\mathcal P$. Let $U$ be a neighborhood of the origin. Since the topology is generated by $\mathcal P$, there is a basic neighborhood
\[
B=\bigcap_{i=1}^{k}\{x\in X\mid p_i(x)<r_i\}
\]
such that $B\subseteq U$. By hypothesis, for each $i$ there exists $M_i>0$ such that $p_i(a)\leq M_i$ for every $a\in A$. Choose $t>0$ such that $M_i<tr_i$ for every $i$. Then $A\subseteq tB\subseteq tU$. Therefore $A$ is topologically bounded.
\end{proof}

\begin{proposition}\label{prop: convergencia y Cauchy por seminormas}
Let $X$ be a vector space, and let $\mathcal P$ be a family, not necessarily countable, of seminorms on $X$. Denote by $\tau_{\mathcal P}$ the topology induced by $\mathcal P$. Then:
\begin{enumerate}[label=(\alph*)]
 \item A sequence $(x_n)_{n\in\mathbb N}$ converges to $x\in X$ in $\tau_{\mathcal P}$ if and only if $p(x_n-x)\to0$ for every $p\in\mathcal P$.
 \item A sequence $(x_n)_{n\in\mathbb N}$ is Cauchy in $\tau_{\mathcal P}$ if and only if $p(x_n-x_m)\to0$ as $n,m\to\infty$, for every $p\in\mathcal P$.
\end{enumerate}
\end{proposition}

\begin{proof}
We first prove the convergence assertion. Suppose that $x_n\to x$ in $\tau_{\mathcal P}$. Let $p\in\mathcal P$ and $\varepsilon>0$. The set
\[
x+\{y\in X\mid p(y)<\varepsilon\}
\]
is a neighborhood of $x$. By convergence of the sequence, there exists $N\in\mathbb N$ such that, if $n\geq N$, then $x_n$ belongs to this neighborhood. Therefore $p(x_n-x)<\varepsilon$ for every $n\geq N$. This proves that $p(x_n-x)\to0$.

Conversely, suppose that $p(x_n-x)\to0$ for every $p\in\mathcal P$. Let $U$ be a neighborhood of $x$. By the definition of the induced topology, there exist $p_1,\dots,p_k\in\mathcal P$ and $r_1,\dots,r_k>0$ such that
\[
x+\bigcap_{i=1}^{k}\{y\in X\mid p_i(y)<r_i\}\subseteq U.
\]
For each $i\in\{1,\dots,k\}$, since $p_i(x_n-x)\to0$, there exists $N_i\in\mathbb N$ such that $p_i(x_n-x)<r_i$ for every $n\geq N_i$. Taking $N:=\displaystyle\max\{N_1,\dots,N_k\}$ gives, for every $n\geq N$,
\[
x_n-x\in\bigcap_{i=1}^{k}\{y\in X\mid p_i(y)<r_i\}.
\]
Therefore $x_n\in U$ for every $n\geq N$. Since $U$ was an arbitrary neighborhood of $x$, we conclude that $x_n\to x$ in $\tau_{\mathcal P}$.

We now prove the assertion about Cauchy sequences. Suppose that $(x_n)_{n\in\mathbb N}$ is Cauchy in $\tau_{\mathcal P}$. Let $p\in\mathcal P$ and $\varepsilon>0$. The set $\{y\in X\mid p(y)<\varepsilon\}$ is a neighborhood of the origin. By the Cauchy property, there exists $N\in\mathbb N$ such that, if $n,m\geq N$, then $x_n-x_m\in\{y\in X\mid p(y)<\varepsilon\}$. Therefore $p(x_n-x_m)<\varepsilon$ for any $n,m\geq N$. This proves that $p(x_n-x_m)\to0$ as $n,m\to\infty$.

Conversely, suppose that, for every $p\in\mathcal P$, we have $p(x_n-x_m)\to0$ as $n,m\to\infty$. Let $U$ be a neighborhood of the origin. By the definition of the induced topology, there exist $p_1,\dots,p_k\in\mathcal P$ and $r_1,\dots,r_k>0$ such that
\[
\bigcap_{i=1}^{k}\{y\in X\mid p_i(y)<r_i\}\subseteq U.
\]
For each $i\in\{1,\dots,k\}$, there exists $N_i\in\mathbb N$ such that $p_i(x_n-x_m)<r_i$ for any $n,m\geq N_i$. Taking $N:=\displaystyle\max\{N_1,\dots,N_k\}$, we have, for any $n,m\geq N$,
\[
x_n-x_m\in
\bigcap_{i=1}^{k}\{y\in X\mid p_i(y)<r_i\}\subseteq U.
\]
Therefore $(x_n)_{n\in\mathbb N}$ is Cauchy in $\tau_{\mathcal P}$.
\end{proof}

\begin{proposition}[Metrization by a countable family of seminorms]\label{prop: metrizacion por seminormas}\index{metrization by a countable family of seminorms}
Let $X$ be a vector space, and let $(p_n)_{n\in\mathbb N_0}$ be a countable family of seminorms separating points. For $x,y\in X$, define
\begin{equation}\label{eq:metrica-frechet-seminormas}
d(x,y):=
\displaystyle\sum_{n=0}^{\infty}
\frac{1}{2^n}
\frac{p_n(x-y)}{1+p_n(x-y)}.
\end{equation}
Then:
\begin{enumerate}[label=(\alph*)]
 \item $d$ is a translation-invariant metric on $X$.
 \item The topology induced by $d$ agrees with the locally convex topology generated by the family $(p_n)_{n\in\mathbb N_0}$.
\end{enumerate}
\end{proposition}

\begin{proof}
Define $\theta\colon [0,\infty)\longrightarrow[0,1)$ by
\[
\theta(t):=\frac{t}{1+t}.
\]
This function is continuous, strictly increasing, satisfies $\theta(0)=0$, and has $\theta(t)=0$ if and only if $t=0$. It is also subadditive. If $s,t\geq0$, then
\[
\theta(s+t)
=
\frac{s+t}{1+s+t}
\leq
\frac{s}{1+s}+\frac{t}{1+t}
=
\theta(s)+\theta(t).
\]

For each $n\in\mathbb N_0$, define
\[
\delta_n(x,y):=\theta(p_n(x-y)),
\qquad x,y\in X.
\]
Since $p_n$ is a seminorm, we have
\[
\delta_n(x+a,y+a)=\delta_n(x,y)
\]
for any $x,y,a\in X$, so $\delta_n$ is translation invariant. Moreover, since $p_n(x-y)=p_n(y-x)$, it is also symmetric.

We show that $\delta_n$ satisfies the triangle inequality. If $x,y,z\in X$, subadditivity of $p_n$ and $\theta$ gives
\[
\delta_n(x,z)=
\theta(p_n(x-z))\leq
\theta(p_n(x-y)+p_n(y-z))\leq
\theta(p_n(x-y))+\theta(p_n(y-z))=
\delta_n(x,y)+\delta_n(y,z).
\]
Thus each $\delta_n$ is a translation-invariant pseudometric. In general, it is not a metric, since a seminorm may vanish on nonzero vectors.

Since $0\leq \delta_n(x,y)<1$ for every $n\in\mathbb N_0$ and any $x,y\in X$, the series
\[
d(x,y):=
\displaystyle\sum_{n=0}^{\infty}\frac{1}{2^n}\delta_n(x,y)
\]
converges, being dominated by the series $\displaystyle\sum_{n=0}^{\infty}2^{-n}$. Symmetry and translation invariance of $d$ follow term by term. Moreover, if $x,y,z\in X$, then
\[
\delta_n(x,z)\leq \delta_n(x,y)+\delta_n(y,z)
\]
for every $n\in\mathbb N_0$. Multiplying by $2^{-n}$ and summing over $n$ gives
\[
d(x,z)\leq d(x,y)+d(y,z).
\]
Thus $d$ satisfies the triangle inequality.

It remains to verify that $d(x,y)=0$ implies $x=y$. If $d(x,y)=0$, all terms of the series are nonnegative, so necessarily $\delta_n(x,y)=0$ for every $n\in\mathbb N_0$. Since $\theta(t)=0$ if and only if $t=0$, it follows that $p_n(x-y)=0$ for every $n\in\mathbb N_0$. The family $(p_n)_{n\in\mathbb N_0}$ separates points, so $x-y=0$, that is, $x=y$. Therefore $d$ is a translation-invariant metric. This proves (a).

We now prove (b). Denote by $\tau$ the locally convex topology generated by the family $(p_n)_{n\in\mathbb N_0}$, and by $\tau_d$ the topology induced by the metric $d$. Since $d$ is translation invariant, the topology $\tau_d$ is also translation invariant. On the other hand, $\tau$ is a topological vector space topology by Theorem~\ref{seminormas topologia} and is therefore also translation invariant. It thus suffices to compare neighborhoods of the origin.

First, we show that every neighborhood of the origin in $\tau_d$ is a neighborhood of the origin in $\tau$. It suffices to prove this for the balls $B_d(0,\varepsilon)$, with $\varepsilon>0$. Choose $N\in\mathbb N_0$ such that
\[
\displaystyle\sum_{n=N+1}^{\infty}\frac{1}{2^n}<\frac{\varepsilon}{2}.
\]
Since $\theta(t)\to0$ as $t\to0^+$, for each $n\in\{0,\dots,N\}$ there exists $\eta_n>0$ such that
\[
\frac{1}{2^n}\theta(\eta_n)<\frac{\varepsilon}{2(N+1)}.
\]
Define
\[
V:=
\bigcap_{n=0}^{N}\{x\in X\mid p_n(x)<\eta_n\}.
\]
Then $V$ is a neighborhood of the origin in $\tau$. If $x\in V$, then $p_n(x)<\eta_n$ for every $n\in\{0,\dots,N\}$, and monotonicity of $\theta$ gives
\[
\begin{aligned}
d(x,0)
&=
\displaystyle\sum_{n=0}^{N}\frac{1}{2^n}\theta(p_n(x))
+
\displaystyle\sum_{n=N+1}^{\infty}\frac{1}{2^n}\theta(p_n(x))\\
&<
\frac{\varepsilon}{2}
+
\displaystyle\sum_{n=N+1}^{\infty}\frac{1}{2^n}
<
\varepsilon.
\end{aligned}
\]
Therefore $V\subseteq B_d(0,\varepsilon)$. Thus every ball centered at the origin for $d$ contains a neighborhood of the origin in $\tau$. Consequently, every neighborhood of the origin in $\tau_d$ is a neighborhood of the origin in $\tau$.

Conversely, we show that every neighborhood of the origin in $\tau$ is a neighborhood of the origin in $\tau_d$. It suffices to consider a basic neighborhood
\[
W:=
\bigcap_{j=1}^{r}\{x\in X\mid p_{n_j}(x)<\varepsilon_j\},
\]
where $n_1,\dots,n_r\in\mathbb N_0$ and $\varepsilon_1,\dots,\varepsilon_r>0$. Define
\[
\rho:=
\min_{1\leq j\leq r}
\frac{1}{2^{n_j}}\theta(\varepsilon_j).
\]
Since each $\varepsilon_j$ is positive and $\theta(t)>0$ for $t>0$, we have $\rho>0$.

We claim that $B_d(0,\rho)\subseteq W$. If $x\in B_d(0,\rho)$, then $d(x,0)<\rho$. For each $j\in\{1,\dots,r\}$, we have
\[
\frac{1}{2^{n_j}}\theta(p_{n_j}(x))
\leq
d(x,0)
<
\rho
\leq
\frac{1}{2^{n_j}}\theta(\varepsilon_j).
\]
Therefore,
\[
\theta(p_{n_j}(x))<\theta(\varepsilon_j).
\]
Since $\theta$ is strictly increasing, $p_{n_j}(x)<\varepsilon_j$. This holds for every $j\in\{1,\dots,r\}$, so $x\in W$. Thus $B_d(0,\rho)\subseteq W$.

We have proved that the neighborhoods of the origin in $\tau$ and $\tau_d$ agree. Since both topologies are translation invariant, their neighborhoods at every point also agree. Consequently, $\tau=\tau_d$.
\end{proof}

\begin{theorem}\label{teo: continuidad mediante seminormas}\index{continuity through seminorms}
Let $X$ be a locally convex topological vector space whose topology is generated by a family of seminorms $\mathcal P$, and let $(Y, \|\cdot\|_Y)$ be a normed space. A linear operator $L\colon X\longrightarrow Y$ is continuous if and only if there exist a constant $C>0$ and finitely many seminorms $p_1, \dots, p_k \in \mathcal P$ such that
\[
\|L(x)\|_Y \leq C \max_{1\leq i\leq k} p_i(x) \qquad \forall x\in X.
\]
In particular, if the family of seminorms is increasing or directed, it suffices to find a single seminorm $p \in \mathcal P$ and a constant $C>0$ such that
\[
\|L(x)\|_Y \leq C p(x) \qquad \forall x\in X.
\]
\end{theorem}

\begin{proof}
First suppose that $L$ is continuous. Then $L^{-1}(B_{d_Y}(0,1))$ is a neighborhood of the origin in $X$. Since the topology of $X$ is generated by $\mathcal P$, there exist $p_1,\dots,p_k\in\mathcal P$ and $\varepsilon>0$ such that
\[
\{x\in X\mid p_i(x)<\varepsilon,\; i=1,\dots,k\}\subseteq L^{-1}(B_{d_Y}(0,1)).
\]
Let $\displaystyle q(x):=\displaystyle\max_{1\leq i\leq k}p_i(x)$. If $x\neq 0$ and $q(x)>0$, then $\displaystyle\frac{\varepsilon}{2q(x)}x$ belongs to the preceding neighborhood, and therefore
\[
\left\|L\left(\frac{\varepsilon}{2q(x)}x\right)\right\|_Y<1.
\]
This gives $\|L(x)\|_Y\leq \displaystyle\frac{2}{\varepsilon}q(x)$. If $q(x)=0$, then $t x$ belongs to the preceding neighborhood for every $t>0$, and linearity yields $L(x)=0$. Thus the inequality holds for every $x\in X$ with $C=\displaystyle\frac{2}{\varepsilon}$.

Conversely, if there exist a constant $C>0$ and seminorms $p_1,\dots,p_k\in\mathcal P$ such that
\[
\|L(x)\|_Y \leq C \max_{1\leq i\leq k} p_i(x),
\]
then, given $\varepsilon>0$, the set
\[
\left\{x\in X\Biggr| p_i(x)<\frac{\varepsilon}{C},\; i=1,\dots,k\right\}
\]
is a neighborhood of the origin mapped into $B_{d_Y}(0,\varepsilon)$. Thus $L$ is continuous at $0$ and, by linearity, at every point.

The final assertion follows because, if the family is increasing or directed, there exists a seminorm $p\in\mathcal P$ that simultaneously dominates $p_1,\dots,p_k$ up to inessential constants.
\end{proof}

The following proposition gives an abstract construction of locally convex topologies from a prescribed local basis.

\begin{proposition}\label{prop: base-discos-topologia-lc}
Let $X$ be a vector space over $\mathbb K$, where $\mathbb K=\mathbb R$ or $\mathbb C$.
Let $\mathcal B$ be a nonempty family of subsets of $X$ satisfying:
\begin{enumerate}[label=(\alph*)]
 \item every $U\in\mathcal B$ is absolutely convex and absorbing;
 \item if $U,V\in\mathcal B$, there exists $W\in\mathcal B$ such that $W\subseteq U\cap V$;
 \item for each $U\in\mathcal B$ and each $\lambda>0$, there exists $W\in\mathcal B$ such that $W\subseteq \lambda U$.
\end{enumerate}
Then there exists a unique topology $\tau$ such that $(X,\tau)$ is a locally convex topological vector space and $\mathcal B$ is a neighborhood basis at the origin.
\end{proposition}

\begin{proof}
For each $U\in\mathcal B$, consider its Minkowski functional
\[
p_U(x):=\inf\{\lambda>0 \mid x\in \lambda U\},
\qquad x\in X.
\]
Since each $U\in\mathcal B$ is absolutely convex and absorbing, Proposition~\ref{prop: minkowski-seminorma} ensures that $p_U$ is a seminorm on $X$.

Define
\[
\mathcal P:=\{p_U\mid U\in\mathcal B\}.
\]
By Theorem~\ref{seminormas topologia}, the family $\mathcal P$ determines a topology $\tau$ making $X$ a locally convex topological vector space. We show that $\mathcal B$ is a local basis at the origin for this topology.

First take $U\in\mathcal B$. If $x\in\{y\in X\mid p_U(y)<1\}$, there exists $\lambda\in(0,1)$ such that $x\in \lambda U$. Since $U$ is balanced, $\lambda U\subseteq U$, and therefore
\[
\{y\in X\mid p_U(y)<1\}\subseteq U.
\]
The set on the left-hand side is a neighborhood of the origin in $\tau$, so $U$ is also a neighborhood of the origin.

It remains to prove that every neighborhood of the origin contains an element of $\mathcal B$. It suffices to consider a basic neighborhood of the form
\[
N=\bigcap_{i=1}^{k}\{x\in X\mid p_{U_i}(x)<r_i\},
\]
where $U_1,\dots,U_k\in\mathcal B$ and $r_1,\dots,r_k>0$. For each $i$, the third hypothesis gives $W_i\in\mathcal B$ such that $W_i\subseteq \displaystyle\frac{r_i}{2}U_i$. If $x\in W_i$, then $p_{U_i}(x)\leq \displaystyle\frac{r_i}{2}<r_i$. Therefore $W_i\subseteq\{x\in X\mid p_{U_i}(x)<r_i\}$. Repeatedly applying the second hypothesis, we find $W\in\mathcal B$ such that
\[
W\subseteq \bigcap_{i=1}^{k}W_i\subseteq N.
\]
Thus every basic neighborhood of the origin contains an element of $\mathcal B$, so $\mathcal B$ is a local basis at the origin for $\tau$.

Finally, uniqueness of $\tau$ follows from Proposition~\ref{prop: topologia-determinada-por-base-local}, since a topological vector space topology is uniquely determined by any local basis at the origin.
\end{proof}

Linear operators were introduced in Definition~\ref{def:b2-espacios-normados-operador-lineal}. If $X$ and $Y$ are topological vector spaces, we retain the notation $\mathcal L(X,Y)$ for the space of continuous linear operators from $X$ to $Y$. In particular, $X':=\mathcal L(X,\mathbb K)$ is the continuous dual of $X$, and its elements are the continuous linear functionals.

\begin{remark}\label{rem: apareamiento dual}
The bilinear map
\[
\langle \cdot,\cdot \rangle_{X',X}\colon X'\times X\longrightarrow \mathbb K,
\qquad
(L,x)\mapsto L(x)
\]
is called the \textit{duality pairing}. The dual space of $X'$, denoted by $X''$, is called the \textit{bidual} of $X$.
\end{remark}

\begin{definition}\label{def: operador acotado entre EVTs}\index{bounded operator!between topological vector spaces}
Let $X$ and $Y$ be topological vector spaces. An operator $L\colon X\longrightarrow Y$ is said to be \textit{bounded} if it maps topologically bounded subsets of $X$ to topologically bounded subsets of $Y$.
\end{definition}

In normed spaces, continuity and boundedness of linear operators are equivalent. In general topological vector spaces, the relationship is more subtle, although useful implications remain. In the metrizable case, continuity can be characterized by sequences.

\begin{theorem}\label{teo: relacion continuidad y acotacion en espacios vectoriales topologicos}\index{relationship between continuity and boundedness in topological vector spaces}
Let $X,Y$ be topological vector spaces, and let $L\colon X\longrightarrow Y$ be a linear operator. Consider the following properties:
\begin{enumerate}[label=(\alph*)]
 \item $L$ is continuous.
 \item $L$ is bounded.
 \item If $\{x_{n}\}_{n\in \mathbb N}$ is a sequence in $X$ and $x_{n}\to 0$, then $\{L(x_{n})\mid n\in \mathbb N\}\subseteq Y$ is a topologically bounded set.
 \item If $\{x_{n}\}_{n\in \mathbb N}$ is a sequence in $X$ and $x_{n}\to 0$, then $L(x_{n})\to 0$.
\end{enumerate}
Then $(a)\Rightarrow (b)\Rightarrow(c)$. If, in addition, $X$ is metrizable, then $(c)\Rightarrow(d)\Rightarrow(a)$ and all four properties are equivalent.
\end{theorem}

\begin{proof}
If $L$ is continuous, it maps topologically bounded sets to topologically bounded sets. Let $A\subseteq X$ be bounded, and let $U$ be a neighborhood of the origin in $Y$. By continuity of $L$ at $0$, there is a neighborhood $V$ of the origin in $X$ such that $L(V)\subseteq U$. Since $A$ is bounded, there exists $t>0$ such that $A\subseteq tV$. By linearity, $L(A)\subseteq tL(V)\subseteq tU$. Therefore $L(A)$ is bounded. This proves $(a)\Rightarrow(b)$.

If $L$ is bounded and $x_n\to 0$ in $X$, the set $\{x_n\mid n\in\mathbb N\}$ is topologically bounded. Its image under $L$ is therefore also topologically bounded. This proves $(b)\Rightarrow(c)$.

Now suppose that $X$ is metrizable and $(c)$ holds. We prove $(d)$ by contradiction. Let $x_n\to 0$ in $X$, and suppose that $L(x_n)$ does not converge to $0$. There is then a balanced neighborhood $U$ of the origin in $Y$ and a subsequence, again denoted by $x_n$, such that $L(x_n)\notin U$ for every $n$.

Choose a decreasing base of open neighborhoods $(V_j)_{j\in\mathbb N}$
of zero in $X$. For every $j$, $(\displaystyle\frac{1}{j})V_j$ is a neighborhood of zero.
Since $x_n\to0$, we can choose indices $n_1<n_2<\cdots$ inductively
with $x_{n_j}\in(\displaystyle\frac{1}{j})V_j$. Thus $z_j:=jx_{n_j}\in V_j$ and $z_j\to0$.
By $(c)$, the set $\{L(z_j)\mid j\in\mathbb N\}$ is topologically bounded.
There is $t>0$ with $jL(x_{n_j})=L(z_j)\in tU$ for every $j$.
For $j>t$, balancedness gives
$L(x_{n_j})\in(\displaystyle\frac{t}{j})U\subseteq U$, contradicting the choice of the
sequence. Hence $L(x_n)\to0$, proving $(c)\Rightarrow(d)$.

Finally, if $X$ is metrizable, the implication $(d)\Rightarrow(a)$ is precisely Proposition~\ref{prop: continuidad secuencial en espacios metrizables}.
\end{proof}

\begin{theorem}[Extension of linear operators]\label{teo: extension continua de operadores}\index{extension of linear operators}
Let $X$ be a metrizable topological vector space, let $D \subseteq X$ be a dense vector subspace, and let $Y$ be a Hausdorff, sequentially complete topological vector space. If $L\colon D\longrightarrow Y$ is a continuous linear operator, where $D$ has the topology induced by $X$, there is at most one continuous linear operator $\widetilde L\colon X\longrightarrow Y$ such that $\widetilde L(x) = L(x)$ for every $x\in D$. Moreover, if, for each $x\in X$ and every sequence $\{x_n\}\subseteq D$ with $x_n\to x$, the sequence $\{L(x_n)\}$ is Cauchy in $Y$, then this extension exists and is unique.
\end{theorem}

\begin{proof}
Uniqueness follows from density of $D$ and the Hausdorff property of $Y$. If $\widetilde L_1$ and $\widetilde L_2$ are two continuous extensions of $L$, they agree on $D$. Since the set where they agree is closed and contains $D$, it is all of $X$.

Now suppose that the additional condition holds. Let $x\in X$. Since $D$ is dense and $X$ is metrizable, there is a sequence $\{x_n\}\subseteq D$ such that $x_n\to x$. By hypothesis, $\{L(x_n)\}$ is Cauchy in $Y$ and, since $Y$ is sequentially complete, it converges to a point of $Y$. Define
\[
\widetilde L(x):=\lim_{n\to\infty}L(x_n).
\]
We show that this definition is independent of the chosen sequence. If $\{y_n\}\subseteq D$ also converges to $x$, the interleaved sequence $z_1=x_1,z_2=y_1,z_3=x_2,z_4=y_2,\dots$ converges to $x$. By hypothesis, $\{L(z_n)\}$ is Cauchy, so its two subsequences $\{L(x_n)\}$ and $\{L(y_n)\}$ have the same limit. Thus $\widetilde L$ is well defined.

Linearity follows by taking sequences in $D$ converging to the corresponding points and using linearity of $L$ before passing to the limit.

Finally, we prove continuity without assuming that $Y$ is first countable. Fix a compatible translation-invariant metric
$d$ on $X$. Let $x_n\to0$, and suppose, for a contradiction, that
$\widetilde L(x_n)\not\to0$. There are a neighborhood $U$ of $0$ in $Y$ and a subsequence, which we do not rename, such that
$\widetilde L(x_n)\notin U$ for every $n$. Choose a neighborhood $V$ of
$0$ such that $V+V\subseteq U$.

For each $n$, the definition of $\widetilde L(x_n)$ gives a sequence of points of $D$ converging to $x_n$ whose images converge to $\widetilde L(x_n)$. Choose a term $d_n\in D$ of this sequence so that
\[
 d(d_n,x_n)<\frac{1}{n},
 \qquad
 \widetilde L(x_n)-L(d_n)\in V.
\]
Then $d_n\to0$ in $X$. Continuity of $L$ at $D$ implies
$L(d_n)\to0$, so $L(d_n)\in V$ for every sufficiently large $n$. For these indices,
\[
 \widetilde L(x_n)
 \in L(d_n)+V
 \subseteq V+V
 \subseteq U,
\]
a contradiction. Therefore $\widetilde L$ is sequentially continuous at the origin. Since $X$ is metrizable,
Proposition~\ref{prop: continuidad secuencial en espacios metrizables} proves that $\widetilde L$ is continuous.
\end{proof}

The following geometric forms of the Hahn--Banach theorem complement the analytical form presented in the section on normed spaces and provide the separation results needed later. Their proofs, obtained by applying the analytical form to the Minkowski functional of a suitable convex open set, can be found in \cite{Rudin1991,Treves1967}.

\begin{theorem}[Hahn--Banach, first geometric form]\label{teo: hahn--banach, primera forma geometrica}\index{Hahn--Banach first geometric form@Hahn--Banach, first geometric form}
Let $X$ be a real locally convex topological vector space, and let
$A,B\subseteq X$ be nonempty convex sets.
Suppose that $\operatorname{Int}(A)\neq\varnothing$ and
$B\cap \operatorname{Int}(A)=\varnothing$.
Then there exists a continuous linear functional $L\colon X\longrightarrow\mathbb R$ such that $\displaystyle \sup_{x\in A}L(x)\le \inf_{y\in B}L(y)$, and moreover $\displaystyle L(a)<\inf_{y\in B}L(y)$ for every $a\in \operatorname{Int}(A)$. In particular, $L\neq 0$.
\end{theorem}

\begin{corollary}\label{cor: dual no trivial}
Let $X$ be a Hausdorff locally convex topological vector space.
Then $X'\neq \{0\}$ if and only if $X$ has a convex neighborhood of the origin strictly contained in $X$.
\end{corollary}

\begin{theorem}[Hahn--Banach, second geometric form]\label{teo: hahn--banach, segunda forma geometrica}\index{hahn banach second geometric form@Hahn--Banach, second geometric form}
Let $X$ be a real locally convex topological vector space, and let $C,K\subseteq X$ be disjoint nonempty convex sets, with $C$ closed and $K$ compact. Then there exist a continuous linear functional $L\colon X\longrightarrow\mathbb R$ and two numbers $a\in \mathbb R$ and $\varepsilon>0$ such that $L(x)\leq a-\varepsilon$ for every $x\in C$ and $L(x)\geq a+\varepsilon$ for every $x\in K$.
\end{theorem}

\begin{theorem}[Mazur]\label{teo: mazur}\index{mazur}
Let $X$ be a topological vector space, and let $M\subseteq X$ be a vector subspace.
Let $K\subseteq X$ be a nonempty open convex subset such that $M\cap K=\varnothing$.

Then there exists a closed hyperplane $N\subseteq X$ such that
\[
M\subseteq N
\qquad\text{and}\qquad
N\cap K=\varnothing .
\]
\end{theorem}

\begin{corollary}[Separation of a subspace and an open convex set]\label{cor: separacion de un subespacio y un abierto convexo, consecuencia de Hahn Banach}\index{separation of a subspace and an open convex set}
Let $X$ be a locally convex Hausdorff topological vector space.
Let $M\subseteq X$ be a vector subspace, and let
$U\subseteq X$ be a nonempty open convex subset such that $M\cap U=\varnothing$.

Then there exists a continuous linear functional
$f\in X'$
such that
\[
f(m)=0 \text{ for every } m\in M,\;\text{and}\;
\operatorname{Re} f(x)>0
\text{ for every } x\in U.
\]
\end{corollary}

\begin{definition}\label{def: encaje continuo EVT}\index{continuous embedding of topological vector spaces}
Let $(X,\tau_X)$ and $(Y,\tau_Y)$ be topological vector spaces. We say that $X$ \textit{is continuously embedded in} $Y$, and write $X\hookrightarrow{Y}$, if $X$ is a subspace of $Y$ and the inclusion $\iota_X\colon X\longrightarrow Y$, given by $\iota_X(x)=x$, is continuous.
\end{definition}

\subsection{Fréchet spaces}\label{ap:espacios-de-Frechet}

Fréchet spaces generalize Banach spaces by allowing the topology to be determined by a countable family of seminorms instead of a single norm. This class arises naturally in the study of spaces of smooth functions and spaces of smooth sections.

\begin{definition}\label{def: espacio de Frechet}\label{def:espacio-frechet-seminormas}\index{frechet space@Fréchet space}
Let $(X,\tau_X)$ be a topological vector space. We say that $X$ is a \textit{Fréchet space} if the following properties hold:
\begin{enumerate}[label=(\alph*)]
 \item $X$ is locally convex.
 \item The topology $\tau_X$ can be induced by a translation-invariant metric; that is, there exists a metric $d\colon X\times X\longrightarrow\mathbb R$ such that
 \[
 d(x+z,y+z)=d(x,y)
 \]
 for any $x,y,z\in X$.
 \item Some, and equivalently every, translation-invariant metric inducing $\tau_X$ is complete.
\end{enumerate}
\end{definition}

\begin{theorem}\label{teo: caracterizacion Frechet por seminormas}\index{frechet space@Fréchet space!characterization by seminorms}
Let $(X,\tau_X)$ be a topological vector space. Then $X$ is a Fréchet space if and only if the following properties hold:
\begin{enumerate}[label=(\alph*)]
 \item $X$ is Hausdorff.
 \item The topology $\tau_X$ is induced by a countable family of seminorms $(p_j)_{j\in\mathbb N_0}$.
 \item $X$ is complete with respect to some, and equivalently every, translation-invariant metric inducing this topology.
\end{enumerate}
\end{theorem}

\begin{proof}
First suppose that $X$ is Fréchet, and fix a complete invariant metric $d$
inducing its topology. For each $j\in\mathbb N_0$,
Proposition~\ref{prop: localmente convexo tiene una base de vecindades absolutamente convexas}
provides an open absolutely convex zero neighborhood $U_j$ contained in
$B_d(0,2^{-j})$. The family $(U_j)_{j\in\mathbb N_0}$ is a local base:
every zero neighborhood contains a ball, hence some $U_j$. By
Theorem~\ref{seminormas topologia}, their Minkowski functionals form a
countable family of seminorms generating the topology. This family
separates points because $X$ is Hausdorff.

Conversely, a countable separating family induces a locally convex
Hausdorff topology and a compatible invariant metric by
Propositions~\ref{prop: Hausdorff seminormas separan puntos} and
\ref{prop: metrizacion por seminormas}. The completeness condition then
gives Definition~\ref{def: espacio de Frechet}.

We specify independence of the metric and the sequential criterion to be
used. For any compatible invariant metric $d$, a sequence is Cauchy if
and only if every zero neighborhood $U$ admits $N$ such that
$x_n-x_m\in U$ for $n,m\geq N$: take the balls $B_d(0,\varepsilon)$ and
use $d(x_n,x_m)=d(x_n-x_m,0)$. This condition depends only on the topology.
By Proposition~\ref{prop: convergencia y Cauchy por seminormas}, it is
equivalent to
\[
 \forall j\in\mathbb N_0\quad\forall\varepsilon>0\quad
 \exists N\in\mathbb N\quad
 n,m\geq N\ \Longrightarrow\ p_j(x_n-x_m)<\varepsilon.
\]
Convergence is likewise tested in every seminorm. Thus completeness is
equivalent to requiring that every sequence satisfying this condition
admit $x\in X$ with $p_j(x_n-x)\to0$ for each $j$; the index $N$ may
depend on $j$ and $\varepsilon$. In particular, two compatible invariant
metrics are complete simultaneously.
\end{proof}

\begin{corollary}[Linear continuity between Fréchet spaces]
\label{prop:metrica-continuidad-frechet}
Let $E,F$ be Fréchet spaces with increasing seminorm families
$(p_j)_{j\in\mathbb N_0}$ and $(q_j)_{j\in\mathbb N_0}$ inducing their
topologies. A linear map $A:E\to F$ is continuous if and only if, for
each $j\in\mathbb N_0$, there are $k_j\in\mathbb N_0$ and $C_j>0$ such that
\[
 q_j(Ax)\leq C_jp_{k_j}(x),\qquad x\in E.
\]
\end{corollary}
\begin{proof}
For each $j$, the vector quotient $F_j:=F/\ker q_j$ is normed by
$\|[y]\|_j=q_j(y)$. This expression is independent of the representative
because $q_j(y+z)=q_j(y)$ for $z\in\ker q_j$.
If $A$ is continuous, its composition with the quotient $F\to F_j$ is
continuous. Theorem~\ref{teo: continuidad mediante seminormas}, applied
to this composition, gives $k_j,C_j$.
Conversely, the bounds make each $q_j\circ A$ continuous. Preimages of
finite intersections of $q_j$-balls are then zero neighborhoods;
linearity proves continuity of $A$.
\end{proof}

Before turning to weak topologies, we recall two general topological constructions. The initial topology will be the natural tool for defining weak topologies, whereas the final topology will arise in the construction of inductive limits.

\subsection{Initial and final topologies}\label{ap:topologias-iniciales-finales}

\begin{definition}[Initial topology]\label{def: topologia inicial}\index{initial topology}
Let $X$ be a set, let $\{(Y_i,\tau_i)\}_{i\in I}$ be a family of topological spaces, and, for each $i\in I$, let $f_i\colon X\longrightarrow Y_i$ be a map. The \textit{initial topology} on $X$ associated with the family $(f_i)_{i\in I}$ is the topology generated by the subbase
\[
\mathcal S_{\mathrm{in}}
:=
\left\{
f_i^{-1}(U)
\ \middle|\
i\in I,\ U\in\tau_i
\right\}.
\]
We denote it by $\tau_{\mathrm{in}}$.
\end{definition}

\begin{proposition}\label{prop: caracterizacion topologia inicial}
The topology $\tau_{\mathrm{in}}$ is the coarsest topology on $X$ for which all the maps $f_i\colon (X,\tau_{\mathrm{in}})\longrightarrow(Y_i,\tau_i)$ are continuous.

Moreover, if $(Z,\tau_Z)$ is a topological space and $g\colon Z\longrightarrow X$ is a map, then $g\colon (Z,\tau_Z)\longrightarrow(X,\tau_{\mathrm{in}})$ is continuous if and only if each composition $f_i\circ g\colon Z\longrightarrow Y_i$ is continuous.
\end{proposition}

\begin{proof}
By the definition of $\tau_{\mathrm{in}}$, if $i\in I$ and $U\subseteq Y_i$ is open, then $f_i^{-1}(U)$ is open in $(X,\tau_{\mathrm{in}})$. Therefore, each map $f_i\colon (X,\tau_{\mathrm{in}})\longrightarrow(Y_i,\tau_i)$ is continuous.

We show that $\tau_{\mathrm{in}}$ is the coarsest topology with this property. Let $\tau$ be a topology on $X$ such that all the maps $f_i\colon (X,\tau)\longrightarrow(Y_i,\tau_i)$ are continuous. Then, if $U$ is open in $Y_i$, the set $f_i^{-1}(U)$ belongs to $\tau$. Since $\tau_{\mathrm{in}}$ is the topology generated by all these sets, it follows that $\tau_{\mathrm{in}}\subseteq\tau$.

We now prove the characterization of continuity. First suppose that $g\colon (Z,\tau_Z)\longrightarrow(X,\tau_{\mathrm{in}})$ is continuous. Since each $f_i$ is continuous, the composition $f_i\circ g\colon Z\longrightarrow Y_i$ is continuous for every $i\in I$.

Conversely, suppose that each composition $f_i\circ g$ is continuous. To prove that $g$ is continuous, it suffices to verify that the preimage under $g$ of each element of the subbase $\mathcal S_{\mathrm{in}}$ is open in $Z$. Let $f_i^{-1}(U)\in\mathcal S_{\mathrm{in}}$, with $U$ open in $Y_i$. Then
\[
g^{-1}\bigl(f_i^{-1}(U)\bigr)
=
(f_i\circ g)^{-1}(U).
\]
The set on the right-hand side is open in $Z$ by the continuity of $f_i\circ g$. Therefore, $g$ is continuous.
\end{proof}

\begin{proposition}\label{prop: convergencia redes topologia inicial}
Let $X$ be a set, let $\{(Y_i,\tau_i)\}_{i\in I}$ be a family of topological spaces, and, for each $i\in I$, let $f_i\colon X\longrightarrow Y_i$ be a map. A net $(x_\alpha)_\alpha$ in $X$ converges to $x\in X$ in the initial topology associated with $(f_i)_{i\in I}$ if and only if
\[
f_i(x_\alpha)\to f_i(x)
\]
in $Y_i$ for every $i\in I$. In particular, a sequence $(x_n)_{n\in\mathbb N}$ converges to $x$ in the initial topology if and only if
\[
f_i(x_n)\to f_i(x)
\]
in $Y_i$ for every $i\in I$.
\end{proposition}

\begin{proof}
If $x_\alpha\to x$ in the initial topology, then, since each $f_i$ is continuous, we have $f_i(x_\alpha)\to f_i(x)$ for every $i\in I$.

Conversely, suppose that $f_i(x_\alpha)\to f_i(x)$ for every $i\in I$. Let $U$ be a neighborhood of $x$ in the initial topology. By the definition of the initial topology, there exists a basic neighborhood $V$ of $x$ such that $x\in V\subseteq U$, where $V$ is a finite intersection of subbasic elements. That is, there exist $i_1,\dots,i_k\in I$ and open sets $U_1\subseteq Y_{i_1},\dots,U_k\subseteq Y_{i_k}$ such that
\[
x\in \bigcap_{j=1}^{k}f_{i_j}^{-1}(U_j)\subseteq U.
\]
Since $f_{i_j}(x_\alpha)\to f_{i_j}(x)$ and $U_j$ is a neighborhood of $f_{i_j}(x)$, for each $j$ there exists $\alpha_j$ such that, if $\alpha\geq \alpha_j$, then $f_{i_j}(x_\alpha)\in U_j$. Since the index set is directed, there exists $\alpha_0$ such that $\alpha_0\geq \alpha_j$ for every $j$. Then, if $\alpha\geq\alpha_0$, we have $x_\alpha\in V\subseteq U$. Therefore, $x_\alpha\to x$.

The statement for sequences is the special case in which the net is indexed by $\mathbb N$.
\end{proof}

\begin{definition}[Final topology]\label{def: topologia final}\index{final topology}
Let $Y$ be a set, let $\{(X_i,\sigma_i)\}_{i\in I}$ be a family of topological spaces, and, for each $i\in I$, let $g_i\colon X_i\longrightarrow Y$ be a map. The \textit{final topology} on $Y$ associated with the family $(g_i)_{i\in I}$ is the collection
\[
\tau_{\mathrm{fin}}
:=
\left\{
U\subseteq Y
\ \middle|\
g_i^{-1}(U)\in\sigma_i
\text{ for every } i\in I
\right\}.
\]
\end{definition}

\begin{proposition}\label{prop: caracterizacion topologia final}
The collection $\tau_{\mathrm{fin}}$ is a topology on $Y$. Moreover, it is the finest topology on $Y$ for which all the maps $g_i\colon (X_i,\sigma_i)\longrightarrow(Y,\tau_{\mathrm{fin}})$ are continuous.

Furthermore, if $(Z,\tau_Z)$ is a topological space and $h\colon Y\longrightarrow Z$ is a map, then $h\colon (Y,\tau_{\mathrm{fin}})\longrightarrow(Z,\tau_Z)$ is continuous if and only if each composition $h\circ g_i\colon X_i\longrightarrow Z$ is continuous.
\end{proposition}

\begin{proof}
First we show that $\tau_{\mathrm{fin}}$ is a topology. Since $g_i^{-1}(\varnothing)=\varnothing$ and $g_i^{-1}(Y)=X_i$ for every $i\in I$, we have $\varnothing,Y\in\tau_{\mathrm{fin}}$.

Now let $(U_\lambda)_{\lambda\in\Lambda}$ be an arbitrary family of elements of $\tau_{\mathrm{fin}}$. For each $i\in I$,
\[
g_i^{-1}\left(\bigcup_{\lambda\in\Lambda}U_\lambda\right)
=
\bigcup_{\lambda\in\Lambda}g_i^{-1}(U_\lambda).
\]
Since each set $g_i^{-1}(U_\lambda)$ is open in $X_i$, the union on the right-hand side is open. Therefore,
\[
\bigcup_{\lambda\in\Lambda}U_\lambda\in\tau_{\mathrm{fin}}.
\]

If $U_1,\dots,U_r\in\tau_{\mathrm{fin}}$, then, for each $i\in I$,
\[
g_i^{-1}\left(\bigcap_{q=1}^{r}U_q\right)
=
\bigcap_{q=1}^{r}g_i^{-1}(U_q).
\]
Since finite intersections of open sets are open, we conclude that
\[
\bigcap_{q=1}^{r}U_q\in\tau_{\mathrm{fin}}.
\]
Thus, $\tau_{\mathrm{fin}}$ is a topology.

By the definition of $\tau_{\mathrm{fin}}$, all the maps $g_i\colon (X_i,\sigma_i)\longrightarrow(Y,\tau_{\mathrm{fin}})$ are continuous.

We show that $\tau_{\mathrm{fin}}$ is the finest topology with this property. Let $\tau$ be a topology on $Y$ such that all the maps $g_i\colon (X_i,\sigma_i)\longrightarrow(Y,\tau)$ are continuous. If $U\in\tau$, then $g_i^{-1}(U)$ is open in $X_i$ for every $i\in I$. By the definition of $\tau_{\mathrm{fin}}$, it follows that $U\in\tau_{\mathrm{fin}}$. Therefore, $\tau\subseteq\tau_{\mathrm{fin}}$.

Finally, we prove the characterization of continuity. If $h\colon (Y,\tau_{\mathrm{fin}})\longrightarrow(Z,\tau_Z)$ is continuous, then each composition $h\circ g_i\colon X_i\longrightarrow Z$ is continuous.

Conversely, suppose that each composition $h\circ g_i$ is continuous. Let $V\subseteq Z$ be an open set. For each $i\in I$,
\[
g_i^{-1}\bigl(h^{-1}(V)\bigr)
=
(h\circ g_i)^{-1}(V),
\]
and this set is open in $X_i$. By the definition of $\tau_{\mathrm{fin}}$, we obtain $h^{-1}(V)\in\tau_{\mathrm{fin}}$. Since this holds for every open set $V\subseteq Z$, we conclude that $h$ is continuous.
\end{proof}

\begin{remark}\label{rem: inicial final debil fuerte}
In some texts, the initial and final topologies associated with a family of maps are called the \textit{weak topology} and the \textit{strong topology}, respectively. In what follows, we will prefer the terms \textit{initial} and \textit{final} to avoid confusion with the weak and weak-$*$ topologies that arise in functional analysis.
\end{remark}

\section{Weak topologies}\label{ap:topologias-debiles}

\subsection{Initial topologies defined by linear operators}

Weak topologies are special cases of the initial topologies introduced in Definition~\ref{def: topologia inicial}. In particular, we shall use the universal property in Proposition~\ref{prop: caracterizacion topologia inicial} and the characterization of convergence in Proposition~\ref{prop: convergencia redes topologia inicial}. When the maps determining the topology are linear, the construction can be described by seminorms.

\begin{proposition}\label{prop: topologia inicial por operadores lineales}
Let $X$ be a vector space over $\mathbb K$, and let $\{Y_i\}_{i\in I}$ be a family of locally convex topological vector spaces. For each $i\in I$, let $T_i\colon X\longrightarrow Y_i$ be a linear operator. Then the initial topology on $X$ associated with $(T_i)_{i\in I}$ is a locally convex vector space topology. If the topology of $Y_i$ is generated by a family of seminorms $\mathcal P_i$, then the initial topology on $X$ is generated by the seminorms
\[
x\longmapsto p(T_i x),
\qquad
i\in I,\quad p\in\mathcal P_i.
\]
\end{proposition}

\begin{proof}
For $i\in I$ and $p\in\mathcal P_i$, define $q_{i,p}\colon X\longrightarrow[0,\infty)$ by $q_{i,p}(x):=p(T_i x)$. Since $T_i$ is linear and $p$ is a seminorm, $q_{i,p}$ is a seminorm on $X$. Let $\tau$ be the topology generated by the family $\{q_{i,p}\mid i\in I,\ p\in\mathcal P_i\}$. By Theorem~\ref{seminormas topologia}, $\tau$ makes $X$ a locally convex topological vector space, and each $T_i\colon (X,\tau)\longrightarrow Y_i$ is continuous.

If $\sigma$ is another topology on $X$ for which all the $T_i$ are continuous, then each seminorm $q_{i,p}=p\circ T_i$ is continuous with respect to $\sigma$. Thus $\tau\subseteq\sigma$, proving that $\tau$ is the initial topology associated with $(T_i)_{i\in I}$.
\end{proof}

\begin{corollary}\label{cor: topologia inicial valores normados}
Let $X$ be a vector space, and let $\mathcal T$ be a family of linear operators $T\colon X\longrightarrow W_T$, where each $W_T$ is a normed space. Then the initial topology on $X$ associated with $\mathcal T$ is generated by the seminorms
\[
x\longmapsto\|Tx\|_{W_T},
\qquad T\in\mathcal T.
\]
In particular, a neighborhood base at the origin consists of finite intersections of sets of the form
\[
\{x\in X\mid \|T_jx\|_{W_{T_j}}<\varepsilon_j,\ j=1,\dots,k\},
\]
where $T_1,\dots,T_k\in\mathcal T$ and $\varepsilon_1,\dots,\varepsilon_k>0$.
\end{corollary}

\begin{proof}
This is a special case of Proposition~\ref{prop: topologia inicial por operadores lineales}, since the topology of a normed space is generated by its norm.
\end{proof}

\begin{corollary}\label{cor: topologia inicial valores finito dimensional}
Under the hypotheses of the preceding corollary, suppose also that each $W_T$ is finite-dimensional. Then the initial topology is also generated by the seminorms
\[
x\longmapsto|\ell(Tx)|,
\qquad T\in\mathcal T,\quad\ell\in W_T'.
\]
\end{corollary}

\begin{proof}
In finite dimensions, the topology of $W_T$ is generated by the linear functionals in $W_T'$. The assertion follows from Proposition~\ref{prop: topologia inicial por operadores lineales} by composing these functionals with the operators $T$.
\end{proof}

\begin{proposition}\label{prop: convergencia topologia inicial normada}
Let $X$ be a vector space, and let $\mathcal T$ be a family of linear operators $T\colon X\longrightarrow W_T$, where each $W_T$ is a normed space. A net $(x_\alpha)_\alpha$ converges to $x\in X$ in the initial topology determined by $\mathcal T$ if and only if $\|T(x_\alpha)-T(x)\|_{W_T}\to0$ for every $T\in\mathcal T$. The same characterization holds for sequences.
\end{proposition}

\begin{proof}
By the characterization of convergence in an initial topology, $x_\alpha\to x$ if and only if $T(x_\alpha)\to T(x)$ in $W_T$ for every $T\in\mathcal T$. Since each $W_T$ is normed, this is equivalent to $\|T(x_\alpha)-T(x)\|_{W_T}\to0$ for every $T\in\mathcal T$.
\end{proof}

\subsection{Weak topologies associated with dualities}

\begin{definition}[Duality]\label{def: dualidad}\index{duality}
Let $X$ and $Y$ be vector spaces over $\mathbb K=\mathbb R$ or $\mathbb C$. A \textit{duality} between $X$ and $Y$ is a bilinear map
\[
\langle\cdot,\cdot\rangle\colon X\times Y\longrightarrow\mathbb K.
\]
We say that the duality separates points of $X$ if, for every $x\in X\setminus\{0\}$, there exists $y\in Y$ such that $\langle x,y\rangle\neq0$. Similarly, it separates points of $Y$ if, for every $y\in Y\setminus\{0\}$, there exists $x\in X$ such that $\langle x,y\rangle\neq0$.
\end{definition}

\begin{definition}[Weak topology associated with a duality]
\label{def: topologia debil dualidad}\index{weak topology associated with a duality}
Let $\langle\cdot,\cdot\rangle\colon X\times Y\longrightarrow\mathbb K$ be a duality. The \textit{weak topology} $\sigma(X,Y)$ on $X$ is the initial topology associated with the linear maps
\[
e_y\colon X\longrightarrow\mathbb K,
\qquad
e_y(x):=\langle x,y\rangle,
\qquad y\in Y.
\]
Equivalently, $\sigma(X,Y)$ is the coarsest topology on $X$ for which all maps $e_y$, with $y\in Y$, are continuous.
\end{definition}

\begin{proposition}[Description by seminorms]
\label{prop: topologia debil dualidad seminormas}\index{description by seminorms}
The topology $\sigma(X,Y)$ is generated by the seminorms
\[
p_y(x):=|\langle x,y\rangle|,
\qquad y\in Y.
\]
In particular, $\sigma(X,Y)$ is locally convex and is Hausdorff if and only if the duality separates points of $X$.
\end{proposition}

\begin{proof}
The first assertion follows from Proposition~\ref{prop: topologia inicial por operadores lineales}, applied to the family $(e_y)_{y\in Y}$. By Proposition~\ref{prop: Hausdorff seminormas separan puntos}, the resulting topology is Hausdorff if and only if $p_y(x)=0$ for every $y\in Y$ implies $x=0$, which is equivalent to the duality separating points of $X$.
\end{proof}

\begin{proposition}\label{prop: convergencia debil dualidad}
Let $\langle\cdot,\cdot\rangle\colon X\times Y\longrightarrow\mathbb K$ be a duality. A net $(x_\alpha)_\alpha$ converges to $x\in X$ with respect to $\sigma(X,Y)$ if and only if
\[
\langle x_\alpha,y\rangle\longrightarrow\langle x,y\rangle
\qquad
\text{for every }y\in Y.
\]
The same characterization holds for sequences.
\end{proposition}

\begin{proof}
The topology $\sigma(X,Y)$ is the initial topology associated with the maps $e_y\colon X\longrightarrow\mathbb K$. The assertion follows from the characterization of convergence in initial topologies.
\end{proof}

If $X$ is a locally convex topological vector space, we may apply the preceding construction to the natural duality between $X$ and its continuous dual $X'$.

\begin{definition}[Weak topology of a locally convex space]
\label{def: topologia debil EVT}\index{weak topology of a locally convex space}
Let $X$ be a locally convex topological vector space. The \textit{weak topology} on $X$ is the topology $\sigma(X,X')$ associated with the duality
\[
\langle x,L\rangle:=L(x),
\qquad
x\in X,\quad L\in X'.
\]
When a sequence $(x_n)_{n\in\mathbb N}$ converges to $x\in X$ in this topology, we write $x_n\rightharpoonup x$.
\end{definition}

\begin{proposition}[Description by seminorms]
\label{prop: topologia debil EVT seminormas}\index{description by seminorms}
The weak topology $\sigma(X,X')$ is generated by the seminorms
\[
p_L(x):=|L(x)|,
\qquad L\in X'.
\]
\end{proposition}

\begin{proof}
This is Proposition~\ref{prop: topologia debil dualidad seminormas} applied to the natural duality between $X$ and $X'$.
\end{proof}

\begin{proposition}\label{prop: caracterizacion convergencia debil}
A net $(x_\alpha)_\alpha$ in $X$ converges weakly to $x\in X$ if and only if $L(x_\alpha)\to L(x)$ for every $L\in X'$. In particular, $x_n\rightharpoonup x$ if and only if
\[
L(x_n)\longrightarrow L(x)
\qquad
\text{for every }L\in X'.
\]
\end{proposition}

\begin{proof}
This is Proposition~\ref{prop: convergencia debil dualidad} applied to the natural duality between $X$ and $X'$.
\end{proof}

The following properties allow weak convergence and weak compactness to be transported by the linear operators arising naturally in applications. The first shows that ordinary linear continuity already includes continuity with respect to weak topologies.

\begin{proposition}[Continuity with respect to weak topologies]
\label{prop:operador-lineal-continuo-topologias-debiles}
\index{weak topology!continuity of linear operators}
Let $X$ and $Y$ be locally convex topological vector spaces, and let
$T\colon X\longrightarrow Y$ be a continuous linear operator. Then
\[
T:\bigl(X,\sigma(X,X')\bigr)
\longrightarrow
\bigl(Y,\sigma(Y,Y')\bigr)
\]
is continuous. In particular, the image under $T$ of a weakly compact set is weakly compact.
\end{proposition}

\begin{proof}
Let $L\in Y'$. Since $T$ is linear and continuous, the composition
$L\circ T\colon X\longrightarrow\mathbb K$ belongs to $X'$. By the definition of
$\sigma(X,X')$, the functional $L\circ T$ is continuous on
$\bigl(X,\sigma(X,X')\bigr)$. Since this holds for every $L\in Y'$ and
$\sigma(Y,Y')$ is the initial topology determined by the functionals in $Y'$,
the universal property of the initial topology implies that
\[
T:\bigl(X,\sigma(X,X')\bigr)
\longrightarrow
\bigl(Y,\sigma(Y,Y')\bigr)
\]
is continuous. The last assertion follows because the continuous image of a
compact set is compact.
\end{proof}

To relate weak convergence obtained through isometric embeddings to the weak topology of the subspace, we shall also use the following consequence of the Hahn--Banach theorem.

\begin{proposition}[Weak topology of a normed subspace]
\label{prop:topologia-debil-subespacio-normado}
\index{weak topology!subspaces}
Let $X$ be a normed space, and let $Y\subseteq X$ be a linear subspace equipped with the induced norm. Then the weak topology $\sigma(Y,Y')$ agrees with the topology induced by $\sigma(X,X')$ on $Y$. In particular, if a net $(y_\alpha)_\alpha\subseteq Y$ and an element $y\in Y$ satisfy
$y_\alpha\rightharpoonup y$ in $X$, then $y_\alpha\rightharpoonup y$ in $Y$, and conversely.
\end{proposition}

\begin{proof}
The topology induced by $\sigma(X,X')$ on $Y$ is generated by the
seminorms
\[
y\longmapsto |L(y)|,
\qquad L\in X'.
\]
Each restriction $L\restriction_Y$ belongs to $Y'$. Conversely, if
$\ell\in Y'$, the Hahn--Banach theorem, in the form of
Corollary~\ref{cor:extension-hahn-banach-preserva-norma}, provides a
functional $L\in X'$ such that $L\restriction_Y=\ell$. Thus
\[
\{L\restriction_Y\mid L\in X'\}=Y',
\]
and the seminorms generating the induced topology are exactly those
generating $\sigma(Y,Y')$. The characterization of convergence follows from
equality of the two topologies.
\end{proof}

The following property describes the weak topology of a finite product. Finiteness is what allows a functional on the product to be recovered from a finite sum of functionals on its components.

\begin{proposition}[Weak topology of a finite product]
\label{prop:topologia-debil-producto-finito}
\index{weak topology!finite products}
Let $X_1,\dots,X_N$ be locally convex topological vector spaces, and let
\[
X:=\prod_{j=1}^N X_j
\]
be equipped with the product topology. Then
\[
\sigma(X,X')
=
\prod_{j=1}^N\sigma(X_j,X_j').
\]
In particular, a net
$x^\alpha=(x_1^\alpha,\dots,x_N^\alpha)$ converges weakly to
$x=(x_1,\dots,x_N)$ in $X$ if and only if
\[
x_j^\alpha\rightharpoonup x_j
\qquad\text{in }X_j
\]
for every $j\in\{1,\dots,N\}$.
\end{proposition}

\begin{proof}
First suppose that $x^\alpha\rightharpoonup x$ in $X$. For each
$j\in\{1,\dots,N\}$, the canonical projection
\[
\pi_j\colon X\longrightarrow X_j
\]
is linear and continuous in the original topologies. By
Proposition~\ref{prop:operador-lineal-continuo-topologias-debiles}, it is also
continuous in the weak topologies. Thus
\[
x_j^\alpha=\pi_j(x^\alpha)
\rightharpoonup\pi_j(x)=x_j
\qquad\text{in }X_j.
\]

Conversely, suppose that
$x_j^\alpha\rightharpoonup x_j$ in $X_j$ for every
$j\in\{1,\dots,N\}$. Let $L\in X'$ and, for each $j$, let
\[
\iota_j\colon X_j\longrightarrow X,
\qquad
\iota_j(z):=(0,\dots,0,z,0,\dots,0),
\]
be the canonical inclusion into the $j$th component. Since $\iota_j$ is linear and
continuous, $L_j:=L\circ\iota_j$ belongs to $X_j'$. Moreover, for every
$(z_1,\dots,z_N)\in X$,
\[
L(z_1,\dots,z_N)
=
\sum_{j=1}^N L_j(z_j).
\]
Consequently,
\[
L(x^\alpha)
=
\sum_{j=1}^N L_j(x_j^\alpha)
\longrightarrow
\sum_{j=1}^N L_j(x_j)
=
L(x),
\]
where passage to the limit is valid because the sum has finitely many
terms. By Proposition~\ref{prop: caracterizacion convergencia debil},
$x^\alpha\rightharpoonup x$ in $X$. Thus the two topologies have exactly
the same convergent nets and therefore agree.
\end{proof}

We also recall a topological characterization of lower semicontinuity that will be useful later.

\begin{proposition}[Characterization by sublevel sets]
\label{prop: semicontinuidad inferior subniveles}\index{characterization by sublevel sets}
Let $(X,\tau)$ be a topological space, and let $F\colon X\longrightarrow(-\infty,+\infty]$. The following assertions are equivalent:
\begin{enumerate}[label=(\alph*)]
\item $F$ is $\tau$-lower semicontinuous.
\item For every $c\in\mathbb R$, the sublevel set
\[
 \{x\in X\mid F(x)\leq c\}
 \]
is $\tau$-closed.
\end{enumerate}
\end{proposition}

\begin{proof}
By the topological characterization of lower semicontinuity, $F$ is $\tau$-lower semicontinuous if and only if $\{x\in X\mid F(x)>c\}$ is $\tau$-open for every $c\in\mathbb R$. Since
\[
\{x\in X\mid F(x)>c\}
=
X\setminus\{x\in X\mid F(x)\leq c\},
\]
this condition is equivalent to all sublevel sets of $F$ being $\tau$-closed.
\end{proof}

\begin{theorem}\label{teo: cerradura convexa fuerte y debil}\index{convex closure!equality of strong and weak closures}
Let $(X,\tau)$ be a locally convex Hausdorff topological vector space, and let $C\subseteq X$ be convex. Then
\[
\overline C^{\,\tau}
=
\overline C^{\,\sigma(X,X')}.
\]
In particular, $C$ is $\tau$-closed if and only if it is weakly closed.
\end{theorem}

\begin{proof}
Since $\sigma(X,X')\subseteq\tau$, we have $\overline C^{\,\tau}\subseteq\overline C^{\,\sigma(X,X')}$. To prove the other inclusion, let $x\notin\overline C^{\,\tau}$. By the Hahn--Banach separation theorem, there exist $L\in X'$ and real numbers $\alpha<\beta$ such that
\[
\operatorname{Re}L(y)\leq\alpha<\beta\leq\operatorname{Re}L(x)
\qquad
\text{for every }y\in\overline C^{\,\tau}.
\]
Then $\{z\in X\mid\operatorname{Re}L(z)>\alpha\}$ is a weak neighborhood of $x$ disjoint from $C$. Thus $x\notin\overline C^{\,\sigma(X,X')}$.
\end{proof}

\begin{corollary}\label{cor: funcional convexo semicontinuo es debilmente semicontinuo}
Let $(X,\tau)$ be a locally convex Hausdorff topological vector space, and let $F\colon X\longrightarrow(-\infty,+\infty]$ be a convex, $\tau$-lower semicontinuous function. Then $F$ is weakly lower semicontinuous.
\end{corollary}

\begin{proof}
By Proposition~\ref{prop: semicontinuidad inferior subniveles}, for every $c\in\mathbb R$ the set $C_c:=\{x\in X\mid F(x)\leq c\}$ is $\tau$-closed. Since $F$ is convex, $C_c$ is convex. By Theorem~\ref{teo: cerradura convexa fuerte y debil}, $C_c$ is weakly closed. Applying Proposition~\ref{prop: semicontinuidad inferior subniveles} again shows that $F$ is weakly lower semicontinuous.
\end{proof}

The uniform boundedness principle was established in Theorem~\ref{teo: banach steinhaus}; we now use it to relate weak boundedness to norm boundedness.

We shall also need the canonical embedding of a normed space into its bidual.

\begin{proposition}\label{prop: inmersion canonica bidual}
Let $X$ be a normed space, and let $J\colon X\longrightarrow X''$ be the linear operator defined by
\[
J(x)(L):=L(x),
\qquad
x\in X,\quad L\in X'.
\]
Then $\|J(x)\|_{X''}=\|x\|_X$ for every $x\in X$. In particular, $J$ is a linear isometry.
\end{proposition}

\begin{proof}
If $\|L\|_{X'}\leq1$, then $|J(x)(L)|=|L(x)|\leq\|x\|_X$, so $\|J(x)\|_{X''}\leq\|x\|_X$. If $x\neq0$, Corollary~\ref{cor:caracterizacion-dual-norma} provides $L\in X'$ such that $\|L\|_{X'}=1$ and $L(x)=\|x\|_X$. Thus $\|J(x)\|_{X''}\geq\|x\|_X$. The case $x=0$ is immediate.
\end{proof}

\begin{theorem}\label{teo: acotacion debil en Banach}\index{weak boundedness in Banach spaces}
Let $X$ be a Banach space, and let $A\subseteq X$. Then $A$ is norm bounded if and only if it is bounded with respect to $\sigma(X,X')$.
\end{theorem}

\begin{proof}
If $A$ is norm bounded, then $\displaystyle \sup_{x\in A}|L(x)|\leq\|L\|_{X'}\displaystyle\sup_{x\in A}\|x\|_X<\infty$ for every $L\in X'$, so $A$ is weakly bounded.

Conversely, if $A$ is weakly bounded, then $\displaystyle \sup_{x\in A}|L(x)|<\infty$ for every $L\in X'$. The family $\{J(x)\mid x\in A\}\subseteq X''$ is therefore pointwise bounded on $X'$. By Theorem~\ref{teo: banach steinhaus},
\[
\sup_{x\in A}\|J(x)\|_{X''}<\infty.
\]
Since $\|J(x)\|_{X''}=\|x\|_X$, we conclude that $A$ is norm bounded.
\end{proof}

\begin{proposition}\label{prop: convergencia debil acotada lsc}
Let $X$ be a Banach space. If $x_n\rightharpoonup x$ in $X$, then $(x_n)_{n\in\mathbb N}$ is bounded and
\[
\|x\|_X\leq\liminf_{n\to\infty}\|x_n\|_X.
\]
\end{proposition}

\begin{proof}
For each $L\in X'$, the sequence $(L(x_n))_{n\in\mathbb N}$ is convergent and hence bounded. Thus $\{x_n\mid n\in\mathbb N\}$ is weakly bounded and, by Theorem~\ref{teo: acotacion debil en Banach}, also norm bounded.

On the other hand, the norm is convex and continuous in the norm topology. By Corollary~\ref{cor: funcional convexo semicontinuo es debilmente semicontinuo}, it is weakly lower semicontinuous. Consequently,
\[
\|x\|_X\leq\liminf_{n\to\infty}\|x_n\|_X.
\]
\end{proof}

\subsection{Weak-$*$ topology and pointwise convergence}

\begin{definition}[Weak-$*$ topology]\label{def: topologia debil estrella}\index{weak star topology@weak-* topology}
Let $X$ be a locally convex topological vector space. The \textit{weak-$*$ topology} on $X'$ is the initial topology $\sigma(X',X)$ associated with the maps
\[
\operatorname{ev}_x\colon X'\longrightarrow\mathbb K,
\qquad
\operatorname{ev}_x(L):=L(x),
\qquad x\in X.
\]
When $L_n\to L$ in this topology, we write $L_n\overset{*}{\rightharpoonup}L$.
\end{definition}

\begin{proposition}[Description by seminorms]
\label{prop: debil estrella seminormas}\index{description by seminorms}
The weak-$*$ topology $\sigma(X',X)$ is generated by the seminorms
\[
p_x(L):=|L(x)|,
\qquad x\in X.
\]
In particular, a net $(L_\alpha)_\alpha$ converges weakly-$*$ to $L$ if and only if
\[
L_\alpha(x)\longrightarrow L(x)
\qquad
\text{for every }x\in X.
\]
The same characterization holds for sequences.
\end{proposition}

\begin{proof}
The description by seminorms follows from Proposition~\ref{prop: topologia debil dualidad seminormas}, applied to the duality between $X'$ and $X$. The characterization of convergence follows from Proposition~\ref{prop: convergencia debil dualidad}.
\end{proof}

When the target space is not the scalar field, the analogous construction is called the topology of pointwise convergence.

\begin{definition}[Topology of pointwise convergence]
\label{def: topologia debil estrella W valuada}\index{topology of pointwise convergence}
Let $X$ be a topological vector space, and let $W$ be a finite-dimensional normed vector space. The \textit{topology of pointwise convergence} on $\mathcal L(X,W)$ is the initial topology associated with the maps
\[
\operatorname{ev}_x\colon \mathcal L(X,W)\longrightarrow W,
\qquad
\operatorname{ev}_x(T):=T(x),
\qquad x\in X.
\]
\end{definition}

\begin{proposition}\label{prop: debil estrella W valuada seminormas}
The topology of pointwise convergence on $\mathcal L(X,W)$ is generated by the seminorms
\[
p_x(T):=\|T(x)\|_W,
\qquad x\in X.
\]
Since $W$ is finite-dimensional, it is also generated by the seminorms
\[
p_{x,\ell}(T):=|\ell(T(x))|,
\qquad x\in X,\quad\ell\in W'.
\]
In particular, a net $(T_\alpha)_\alpha$ converges to $T$ in this topology if and only if
\[
T_\alpha(x)\longrightarrow T(x)
\qquad
\text{in }W
\]
for every $x\in X$. The same characterization holds for sequences.
\end{proposition}

\begin{proof}
The first description follows from Corollary~\ref{cor: topologia inicial valores normados}, applied to the maps $\operatorname{ev}_x$. Since $W$ is finite-dimensional, its topology is generated by the seminorms $w\mapsto|\ell(w)|$, with $\ell\in W'$, giving the second description. The characterization of convergence follows from Proposition~\ref{prop: convergencia topologia inicial normada}.
\end{proof}

\begin{remark}\label{obs: top debil estrella distribuciones W valuadas abstracta}
This construction arises naturally in spaces of distributions with values in a finite-dimensional vector space. If $\mathcal E$ is a test space, the space $\mathcal L(\mathcal E,W)$ is equipped with the topology of pointwise convergence, so $F_\alpha\to F$ if and only if
\[
F_\alpha(\varphi)\longrightarrow F(\varphi)
\qquad
\text{in }W
\]
for every $\varphi\in\mathcal E$. This is precisely weak convergence of $W$-valued distributions.
\end{remark}

\subsection{Weak compactness and reflexivity}

The Banach--Alaoglu, Eberlein--Šmulian, and Kakutani theorems will be used in the forms stated below; their proofs can be found in \cite{brezis_functional_2011,Rudin1991}. In this subsection we prove the weak and weak-\(*\) metrizability criteria, together with the consequences of compactness and reflexivity used later.

\begin{definition}\label{def: polar}\index{polar}
Let $X$ be a topological vector space, and let $A\subseteq X$. The \textit{polar} of $A$ is the set
\[
A^\circ
:=
\{L\in X'\mid |L(x)|\leq1\text{ for every }x\in A\}.
\]
\end{definition}

\begin{theorem}[Banach--Alaoglu]\label{teo: banach--alaoglu}\index{Banach--Alaoglu}
If $U$ is a neighborhood of the origin in a locally convex topological vector space $X$, then $U^\circ$ is compact in the weak-$*$ topology $\sigma(X',X)$.
\end{theorem}

\begin{corollary}\label{cor: banach alaoglu bola dual}
If $X$ is a normed space, then the closed unit ball of $X'$ is weakly-$*$ compact.
\end{corollary}

\begin{proof}
The closed unit ball of $X'$ is the polar of the open unit ball of $X$. The assertion follows from Theorem~\ref{teo: banach--alaoglu}.
\end{proof}

\begin{theorem}\label{teo: K es metrizable topologia debil * estrella}\index{weak-* topology!metrizability on bounded compact sets}
\label{teo: metrizabilidad debil estrella en acotados}
Let $X$ be a separable normed space, and let $K\subseteq X'$ be a bounded weakly-$*$ compact set. Then the weak-$*$ topology induced on $K$ is metrizable. In particular, the closed unit ball of $X'$ is weakly-$*$ compact and metrizable.
\end{theorem}

\begin{proof}
Let $\{x_j\}_{j\in\mathbb N}$ be a countable dense subset of $X$. For $L_1,L_2\in K$, define
\[
d(L_1,L_2)
:=
\displaystyle\sum_{j=1}^{\infty}
2^{-j}
\frac{|(L_1-L_2)(x_j)|}{1+|(L_1-L_2)(x_j)|}.
\]
If $d(L_1,L_2)=0$, then $L_1-L_2$ vanishes on the dense set $\{x_j\mid j\in\mathbb N\}$ and, by continuity, on all of $X$. The triangle inequality follows because $t\mapsto \displaystyle\frac{t}{1+t}$ is
increasing and subadditive: for $a,b\geq0$,
\[
 \frac{a+b}{1+a+b}
 =\frac{a}{1+a+b}+\frac{b}{1+a+b}
 \leq\frac{a}{1+a}+\frac{b}{1+b}.
\]
Thus $d$ is a metric.

Each summand is continuous in the weak-$*$ topology, and the series converges uniformly because it is dominated by $\displaystyle\sum_{j=1}^{\infty}2^{-j}$. Consequently, the identity from $K$ with the weak-$*$ topology to $K$ with the topology induced by $d$ is continuous. Since the first space is compact and the second Hausdorff, this identity is a homeomorphism. Thus $d$ metrizes the weak-$*$ topology on $K$.
\end{proof}

\begin{corollary}\label{cor: subsucesion debil estrella}
Let $X$ be a separable Banach space. Every bounded sequence in $X'$ admits a weakly-$*$ convergent subsequence.
\end{corollary}

\begin{proof}
Every bounded sequence lies in a closed ball of $X'$. By Theorems~\ref{teo: banach--alaoglu} and~\ref{teo: metrizabilidad debil estrella en acotados}, this ball is weakly-$*$ compact and metrizable. The conclusion follows from sequential compactness of compact metric spaces.
\end{proof}

\begin{theorem}\label{teo: bola dual debil estrella metrizable iff separable}\index{weak-* topology!metrizability of the dual ball}
Let $X$ be a Banach space. The closed unit ball of $X'$, equipped with the weak-$*$ topology, is metrizable if and only if $X$ is separable.
\end{theorem}

\begin{proof}
If $X$ is separable, Corollary~\ref{cor: banach alaoglu bola dual} and
Theorem~\ref{teo: metrizabilidad debil estrella en acotados} give the
conclusion. Conversely, let $B_{X'}$ be the closed dual ball, and suppose
its weak-$*$ topology is metrizable. Zero has a countable relative
neighborhood base $(V_j)_{j\in\mathbb N}$. For every $j$ there are a
finite set $F_j\subseteq X$ and $\varepsilon_j>0$ such that
\[
 W_j:=\{L\in B_{X'}\mid |L(x)|<\varepsilon_j\text{ for every }x\in F_j\}
 \subseteq V_j.
\]
Set $Y:=\overline{\operatorname{span}\bigcup_jF_j}^{\,\|\cdot\|_X}$.
If $Y\neq X$, Hahn--Banach gives a norm-one functional $L\in X'$ that
vanishes on $Y$. Then $L\in W_j\subseteq V_j$ for every $j$.
The weak-$*$ topology is Hausdorff and $(V_j)$ is a base at zero,
so $\bigcap_jV_j=\{0\}$, contradicting $\|L\|=1$.
Thus $Y=X$. The span of the countable set $\bigcup_jF_j$ with rational
coefficients, or complex rational coefficients, is countable and dense
in $X$. Therefore $X$ is separable.
\end{proof}

\begin{proposition}\label{prop: debil estrella acotada lsc}
Let $X$ be a Banach space. If $L_n\overset{*}{\rightharpoonup}L$ in $X'$, then $(L_n)_{n\in\mathbb N}$ is bounded and
\[
\|L\|_{X'}
\leq
\liminf_{n\to\infty}\|L_n\|_{X'}.
\]
\end{proposition}

\begin{proof}
For each $x\in X$, the sequence $(L_n(x))_{n\in\mathbb N}$ is convergent and hence bounded. Theorem~\ref{teo: banach steinhaus} implies that $\displaystyle \sup_{n\in\mathbb N}\|L_n\|_{X'}<\infty$.

By Theorem~\ref{teo: banach--alaoglu}, closed balls in $X'$ are weakly-$*$ compact and, since the weak-$*$ topology is Hausdorff, weakly-$*$ closed. These sets are precisely the sublevel sets of the norm. By Proposition~\ref{prop: semicontinuidad inferior subniveles}, the norm on $X'$ is weakly-$*$ lower semicontinuous, giving the inequality.
\end{proof}

\begin{proposition}\label{prop: dual separable implica espacio separable}
Let $X$ be a Banach space. If $X'$ is separable, then $X$ is separable.
\end{proposition}
\begin{proof}
If \(X'=\{0\}\), the dual characterization of the norm, Corollary~\ref{cor:caracterizacion-dual-norma}, implies that \(X=\{0\}\). Suppose \(X'\neq\{0\}\). Since \(X'\) is separable, its unit sphere contains a dense sequence \((\varphi_n)_{n\in\mathbb N}\). For each \(n\), by the definition of the dual norm, we may choose \(x_n\in X\), \(\|x_n\|_X=1\), such that
\[
|\varphi_n(x_n)|>\frac{1}{2}.
\]
Let \(Y:=\overline{\operatorname{span}}\{x_n\mid n\in\mathbb N\}\). If \(Y\neq X\), Hahn--Banach provides \(\varphi\in X'\) such that \(\|\varphi\|=1\) and \(\varphi\restriction_Y=0\). Choose \(n\) with \(\|\varphi_n-\varphi\|<\displaystyle\frac{1}{2}\). Since \(x_n\in Y\),
\[
|\varphi_n(x_n)|
=|(\varphi_n-\varphi)(x_n)|
\leq\|\varphi_n-\varphi\|\,\|x_n\|
<\frac{1}{2},
\]
which contradicts the choice of \(x_n\). Thus \(Y=X\). The subspace spanned by \(\{x_n\}\) using rational coefficients---complex rational coefficients if \(\mathbb K=\mathbb C\)---is countable and dense in \(X\), so \(X\) is separable.
\end{proof}

The converse is false. For example, $L^{1}(\mathbb R^n)$ is separable, but its dual $L^{\infty}(\mathbb R^n)$ is not.

\begin{theorem}\label{teo: bola debil metrizable dual separable}\index{unit ball!weak metrizability}
Let $X$ be a Banach space whose dual $X'$ is separable. Then the closed unit ball of $X$, equipped with the weak topology $\sigma(X,X')$, is metrizable.
\end{theorem}

\begin{proof}
Let $(\varphi_j)_{j\in\mathbb N}$ be norm dense in $X'$.
On the closed ball $B_X$ define
\[
 d(x,y):=\sum_{j=1}^{\infty}2^{-j}
       \frac{|\varphi_j(x-y)|}{1+|\varphi_j(x-y)|}.
\]
The function $f(t)=\displaystyle\frac{t}{1+t}$ is increasing and subadditive on $[0,+\infty)$:
\[
 f(a+b)=\frac{a}{1+a+b}+\frac{b}{1+a+b}
 \leq\frac{a}{1+a}+\frac{b}{1+b}=f(a)+f(b).
\]
This proves the triangle inequality for $d$. If $d(x,y)=0$, all the
$\varphi_j$ vanish on $x-y$. Density gives the same conclusion for
every functional in $X'$, and Hahn--Banach implies $x=y$.
Thus $d$ is a metric.

Every summand is weakly continuous and the series converges uniformly,
so the identity from the weak topology on $B_X$ to the topology of $d$
is continuous. For the converse, fix $x\in B_X$ and a basic weak
neighborhood
\[
 V:=\{y\in B_X\mid |\psi_i(y-x)|<\varepsilon,
                       \ i=1,\ldots,m\},
 \qquad \varepsilon>0.
\]
Choose indices $j_i$ with $\|\psi_i-\varphi_{j_i}\|_{X'}<\displaystyle\frac{\varepsilon}{4}$.
Since $\|y-x\|_X\leq2$ on $B_X$, we have
\[
 |\psi_i(y-x)|<\frac{\varepsilon}{2}+|\varphi_{j_i}(y-x)|.
\]
If $d(x,y)<\min_i 2^{-j_i}f(\displaystyle\frac{\varepsilon}{2})$, each corresponding term
satisfies $|\varphi_{j_i}(y-x)|<\displaystyle\frac{\varepsilon}{2}$.
This metric ball is therefore contained in $V$, proving that the
two topologies agree.
\end{proof}

\begin{definition}\label{def: secuencialmente debilmente compacto}\index{weakly sequentially compact}
Let $X$ be a locally convex topological vector space. A set $K\subseteq X$ is \textit{weakly sequentially compact} if every sequence in $K$ admits a subsequence converging weakly to a point of $K$.
\end{definition}

\begin{theorem}[Eberlein--Šmulian]\label{teo: eberlein smulian}\index{Eberlein--Šmulian}
Let $X$ be a Banach space, and let $K\subseteq X$. Then $K$ is weakly compact if and only if it is weakly sequentially compact.
\end{theorem}

\begin{definition}\label{def: reflexivo}\index{reflexive}
A Banach space $X$ is \textit{reflexive} if the canonical embedding $J\colon X\longrightarrow X''$ is surjective.
\end{definition}

In this case, the canonical embedding identifies $X$ isometrically with its bidual $X''$.

\begin{theorem}[Kakutani]\label{teo: kakutani}\index{Kakutani}
A Banach space is reflexive if and only if its closed unit ball is weakly compact.
\end{theorem}

\begin{corollary}\label{cor: sucesion acotada reflexivo tiene subsucesion debil}
Let $X$ be a reflexive Banach space. Every bounded sequence in $X$ admits a weakly convergent subsequence.
\end{corollary}

\begin{proof}
Every bounded sequence lies in a closed ball. By Theorem~\ref{teo: kakutani}, this ball is weakly compact and, by Theorem~\ref{teo: eberlein smulian}, weakly sequentially compact.
\end{proof}

\begin{proposition}\label{prop: reflexivo iff dual reflexivo}
A Banach space $X$ is reflexive if and only if $X'$ is reflexive.
\end{proposition}
\begin{proof}
First suppose that \(X\) is reflexive, and let \(\Phi\in X'''\). Define \(\varphi\in X'\) by
\[
\varphi(x):=\Phi(J_Xx),
\qquad x\in X,
\]
where \(J_X\colon X\to X''\) is the canonical embedding. Given \(x''\in X''\), reflexivity provides \(x\in X\) such that \(x''=J_Xx\). Then
\[
(J_{X'}\varphi)(x'')
=x''(\varphi)
=\varphi(x)
=\Phi(J_Xx)
=\Phi(x'').
\]
Thus \(J_{X'}\) is surjective and \(X'\) is reflexive.

Conversely, suppose that \(X'\) is reflexive. The image \(J_X(X)\) is a closed subspace of \(X''\), since \(J_X\) is an isometry and \(X\) is complete. If \(J_X(X)\neq X''\), the Hahn--Banach theorem provides a nonzero functional \(\Phi\in X'''\) vanishing on \(J_X(X)\). By reflexivity of \(X'\), there exists \(\varphi\in X'\) such that \(\Phi=J_{X'}\varphi\). For each \(x\in X\),
\[
0=\Phi(J_Xx)=(J_Xx)(\varphi)=\varphi(x),
\]
so \(\varphi=0\) and hence \(\Phi=0\), a contradiction. Consequently, \(J_X\) is surjective and \(X\) is reflexive.
\end{proof}

\section{Topological tensor products and nuclear spaces}
\label{ap:productos-tensoriales-topologicos-nucleares}

In this section we develop the part of topological tensor product theory needed for the Schwartz kernel theorem. The distinction between the projective, injective, and inductive topologies is essential: on the nuclear Fréchet spaces arising when supports are fixed, the projective and injective topologies agree; on the global $LF$ spaces, however, the product representing separate continuity is the inductive one. The exposition follows Chapters~43, 50, and 51 of \cite{Treves1967}. We include proofs of the elementary results used later and explicitly identify the single hereditary theorem on nuclearity taken from that reference.

\subsection{The projective, injective, and inductive topologies}

Let $X$ and $Y$ be vector spaces over $\mathbb K$. Recall that the algebraic tensor product $X\otimes Y$, together with the canonical bilinear map $(x,y)\mapsto x\otimes y$, is characterized by the universal property in Definition~\ref{def:nociones-fundamentales-propiedad-universal-del-producto-tensorial}. When $X$ and $Y$ are locally convex, there are several natural topologies on this same algebraic vector space.

\begin{definition}[Projective tensor topology]
\label{def:topologia-tensorial-proyectiva}
\index{projective tensor product}
Let $X$ and $Y$ be locally convex spaces. If $p$ is a continuous seminorm on $X$ and $q$ is a continuous seminorm on $Y$, define, for $\mathbf{u}\in X\otimes Y$,
\[
(p\otimes_{\pi}q)(\mathbf{u})
:=
\inf
\left\{
\displaystyle\sum_{j=1}^{r}p(x_j)q(y_j)
\ \middle|\
\mathbf{u}=\displaystyle\sum_{j=1}^{r}x_j\otimes y_j
\right\}.
\]
The locally convex topology generated by these seminorms is called the \textit{projective tensor topology}; the resulting space is denoted by $X\otimes_{\pi}Y$ and its Hausdorff completion by $X\widehat\otimes_{\pi}Y$.
\end{definition}

The infimum is taken over a nonempty set, since every element of the algebraic tensor product is a finite sum of elementary tensors. If one of the spaces is not Hausdorff, completion will always mean the completion of the quotient by the closure of $\{0\}$; all spaces in the applications will be Hausdorff.

\begin{proposition}[Universal property of the projective product]
\label{prop:propiedad-universal-producto-proyectivo}
For every pair of continuous seminorms $p$ and $q$, the map $p\otimes_{\pi}q$ is a seminorm and satisfies
\[
(p\otimes_{\pi}q)(x\otimes y)=p(x)q(y).
\]
Moreover, the projective topology is the finest locally convex topology on $X\otimes Y$ for which the canonical bilinear map $X\times Y\to X\otimes Y$ is continuous.

More generally, if $Z$ is a locally convex space and $b\colon X\times Y\longrightarrow Z$ is bilinear, then $b$ is continuous if and only if its algebraic linearization
\[
\widetilde b\colon X\otimes Y\longrightarrow Z,
\qquad
\widetilde b(x\otimes y):=b(x,y),
\]
is continuous on $X\otimes_{\pi}Y$. If $Z$ is complete and Hausdorff, $\widetilde b$ extends uniquely to a continuous operator $X\widehat\otimes_{\pi}Y\to Z$.
\end{proposition}

\begin{proof}
Homogeneity of $p\otimes_\pi q$ follows by multiplying one factor
of each representation by the corresponding scalar. Concatenating
representations of $u$ and $v$ and taking infima proves the triangle
inequality. For an elementary tensor, the definition gives
$(p\otimes_\pi q)(x\otimes y)\leq p(x)q(y)$. If either factor is
zero, equality is already proved. Otherwise, Hahn--Banach on the
normed spaces $X/\ker p$ and $Y/\ker q$ gives functionals $\lambda$
and $\mu$ of norm one with $\lambda([x])=p(x)$ and $\mu([y])=q(y)$.
For every representation $x\otimes y=\sum_jx_j\otimes y_j$,
\[
 p(x)q(y)=\left|\sum_j\lambda([x_j])\mu([y_j])\right|
 \leq\sum_jp(x_j)q(y_j).
\]
Taking the infimum gives the reverse inequality.

Let $\rho:=p\otimes_\pi q$. The elementary identity and bilinearity give
\[
 \rho(x\otimes y-x_0\otimes y_0)
 \leq p(x-x_0)q(y_0)+p(x_0)q(y-y_0)
       +p(x-x_0)q(y-y_0).
\]
Continuity of $p$ and $q$ therefore makes the canonical map continuous
at every $(x_0,y_0)$ for the projective topology.

Let $\tau$ be a locally convex topology making that map continuous,
and let $r$ be a $\tau$-continuous seminorm. There are absolutely
convex neighborhoods $U$ and $V$ of zero such that
$r(x\otimes y)\leq1$ for $x\in U$, $y\in V$.
Let $p_U$ and $q_V$ be their Minkowski functionals. If
$\alpha>p_U(x)$ and $\beta>q_V(y)$, then $\displaystyle\frac{x}{\alpha}\in U$,
$\displaystyle\frac{y}{\beta}\in V$, and
\[
 r(x\otimes y)\leq\alpha\beta.
\]
Letting $\alpha\downarrow p_U(x)$ and $\beta\downarrow q_V(y)$
gives $r(x\otimes y)\leq p_U(x)q_V(y)$, including the cases
in which either seminorm vanishes. For $u=\sum_jx_j\otimes y_j$,
the triangle inequality and the infimum over representations give
$r(u)\leq(p_U\otimes_\pi q_V)(u)$. This proves maximality of
the projective topology.

If $b$ is continuous, the same argument applied to each continuous
seminorm $r$ on $Z$ gives continuous seminorms $p$ and $q$ with
\[
 r(b(x,y))\leq p(x)q(y),\qquad
 r(\widetilde b(u))\leq(p\otimes_\pi q)(u).
\]
Hence $\widetilde b$ is continuous. The converse follows by composing
$\widetilde b$ with the canonical bilinear map.

Now assume that $Z$ is complete and Hausdorff. If
$N:=\overline{\{0\}}\subseteq X\otimes_\pi Y$, every continuous
seminorm vanishes on $N$. The last estimate and separation by
seminorms in $Z$ imply $\widetilde b(N)=\{0\}$.
Thus $\widetilde b$ induces a continuous operator on the Hausdorff
quotient $(X\otimes_\pi Y)/N$. For an element of its completion,
choose a net in the quotient converging to that element. This net
is Cauchy, and so is its image, since every continuous linear
operator is uniformly continuous. Completeness of $Z$ defines
the extension as the limit of those images. If another net is
chosen, the difference converges to zero and so does its image;
the limit is therefore independent of the choice. Linearity and
the seminorm bounds pass to the limit, making the extension linear
and continuous. Density of the quotient and the Hausdorff property
of $Z$ prove uniqueness.
\end{proof}

\begin{proposition}[Projective product and finite products]
\label{prop:tensor-proyectivo-productos-finitos}
Let $X_1,\dots,X_A$ and $Y_1,\dots,Y_B$ be locally convex Hausdorff spaces. If $\operatorname{pr}_a$ and $\operatorname{pr}_b$ denote the corresponding coordinate projections, the operator
\[
\Theta:
\left(\prod_{a=1}^{A}X_a\right)
\otimes_{\pi}
\left(\prod_{b=1}^{B}Y_b\right)
\longrightarrow
\prod_{a=1}^{A}\prod_{b=1}^{B}
\left(X_a\otimes_{\pi}Y_b\right),
\qquad
\Theta(u):=
\bigl((\operatorname{pr}_a\otimes\operatorname{pr}_b)u\bigr)_{a,b},
\]
is a topological isomorphism. Consequently, it induces a topological isomorphism between the completions
\[
\left(\prod_{a=1}^{A}X_a\right)
\widehat\otimes_{\pi}
\left(\prod_{b=1}^{B}Y_b\right)
\cong
\prod_{a=1}^{A}\prod_{b=1}^{B}
\left(X_a\widehat\otimes_{\pi}Y_b\right).
\]
\end{proposition}

\begin{proof}
Let $\displaystyle \iota_a\colon X_a\longrightarrow\prod_cX_c$ and $\displaystyle \jmath_b\colon Y_b\longrightarrow\prod_dY_d$ be the coordinate inclusions. For each pair $(a,b)$, the bilinear map
\[
(x,y)\longmapsto
\operatorname{pr}_a(x)\otimes\operatorname{pr}_b(y)
\]
is continuous with values in $X_a\otimes_{\pi}Y_b$: it is the composition of the continuous projections with the canonical bilinear map. Proposition~\ref{prop:propiedad-universal-producto-proyectivo} shows that $\operatorname{pr}_a\otimes\operatorname{pr}_b$ is continuous; since there are only finitely many pairs, $\Theta$ is continuous as well.

Algebraically, the inverse of $\Theta$ is given by
\[
\Psi\bigl((u_{a,b})_{a,b}\bigr)
:=
\displaystyle\sum_{a=1}^{A}\displaystyle\sum_{b=1}^{B}
(\iota_a\otimes\jmath_b)u_{a,b}.
\]
The same argument proves that each $\iota_a\otimes\jmath_b$ is continuous and, since the sum is finite, $\Psi$ is continuous. For $x=(x_a)_a$ and $y=(y_b)_b$, we have
\[
\Psi\Theta(x\otimes y)
=
\displaystyle\sum_{\substack{1\leq a\leq A\\1\leq b\leq B}}\iota_a(x_a)\otimes\jmath_b(y_b)
=
\left(\displaystyle\sum_{a=1}^{A}\iota_a(x_a)\right)
\otimes
\left(\displaystyle\sum_{b=1}^{B}\jmath_b(y_b)\right)
=x\otimes y,
\]
while the identity $\Theta\Psi=\operatorname{Id}$ is checked in each coordinate. Since elementary tensors span the algebraic tensor product, $\Theta$ and $\Psi$ are inverses.

Every topological isomorphism between locally convex Hausdorff spaces extends, together with its inverse, to their completions. Finally, a net is Cauchy in a finite product if and only if each component is Cauchy; thus the completion of $\displaystyle \prod_{a,b}(X_a\otimes_{\pi}Y_b)$ is canonically identified with $\displaystyle \prod_{a,b}(X_a\widehat\otimes_{\pi}Y_b)$. This proves the second assertion.
\end{proof}

To define the injective topology, we need the notion of equicontinuity.

\begin{definition}[Equicontinuous family]
\label{def:familia-equicontinua-dual}
\index{equicontinuity!in the continuous dual}
Let $X$ be a locally convex space. A subset $A\subseteq X'$ is \textit{equicontinuous} if there exists a neighborhood $U$ of the origin in $X$ such that
\[
|\lambda(x)|\leq1
\qquad
\text{for every }\lambda\in A\text{ and every }x\in U.
\]
Equivalently, the map $\displaystyle p_A(x):=\displaystyle\sup_{\lambda\in A}|\lambda(x)|$ is a continuous seminorm on $X$.
\end{definition}

\begin{proof}[Justification of the equivalence]
If $A$ satisfies the first condition, then $p_A\leq p_U$, where $p_U$ is the Minkowski functional of an absolutely convex neighborhood contained in $U$; thus $p_A$ is continuous. Conversely, if $p_A$ is continuous, the set $\{x\in X\mid p_A(x)<1\}$ is a neighborhood of the origin and $|\lambda(x)|<1$ on it for every $\lambda\in A$.
\end{proof}

Every tensor $\mathbf{u}=\displaystyle\sum_{j=1}^{r}x_j\otimes y_j$ determines a bilinear form on $X'\times Y'$ by
\[
\mathbf{u}(\lambda,\mu)
:=
\displaystyle\sum_{j=1}^{r}\lambda(x_j)\mu(y_j).
\]
The algebraic universal property shows that this expression is independent of the chosen representation.

\begin{definition}[Injective tensor topology]
\label{def:topologia-tensorial-inyectiva}
\index{injective tensor product}
Let $X$ and $Y$ be locally convex Hausdorff spaces. For equicontinuous subsets $A\subseteq X'$ and $B\subseteq Y'$, define
\[
\varepsilon_{A,B}(u)
:=
\sup_{\lambda\in A,\,\mu\in B}|u(\lambda,\mu)|.
\]
The locally convex topology generated by the seminorms $\varepsilon_{A,B}$ is called the \textit{injective tensor topology}. Denote the resulting space by $X\otimes_{\varepsilon}Y$ and its completion by $X\widehat\otimes_{\varepsilon}Y$.
\end{definition}

\begin{proposition}[Comparison between $\varepsilon$ and $\pi$]
\label{prop:epsilon-menos-fina-pi}
The injective topology is coarser than the projective topology. More precisely, if $A\subseteq X'$ and $B\subseteq Y'$ are equicontinuous, there exist continuous seminorms $p$ on $X$ and $q$ on $Y$ such that
\[
\varepsilon_{A,B}(\mathbf{u})
\leq
(p\otimes_{\pi}q)(\mathbf{u})
\qquad
\text{for every }\mathbf{u}\in X\otimes Y.
\]
In particular, the algebraic identity induces continuous operators
\[
X\otimes_{\pi}Y\longrightarrow X\otimes_{\varepsilon}Y,
\qquad
X\widehat\otimes_{\pi}Y\longrightarrow X\widehat\otimes_{\varepsilon}Y.
\]
\end{proposition}

\begin{proof}
Equicontinuity gives continuous seminorms $p$ and $q$ such that
$|\lambda(x)|\leq p(x)$ and $|\mu(y)|\leq q(y)$ for
$\lambda\in A$ and $\mu\in B$. If $u=\sum_{j=1}^Nx_j\otimes y_j$,
\[
 |u(\lambda,\mu)|\leq\sum_{j=1}^Np(x_j)q(y_j).
\]
Taking the supremum over $A\times B$ and then the infimum over
representations of $u$ gives the stated inequality.

These topologies are Hausdorff when $X$ and $Y$ are Hausdorff.
Indeed, if $u\neq0$, write $u=\sum_{j=1}^Nx_j\otimes y_j$
with the $x_j$ linearly independent and some $y_j\neq0$.
The coordinate functionals on $\operatorname{span}\{x_1,\ldots,x_N\}$
are continuous, and Hahn--Banach extends them to functionals in $X'$.
Choose $\lambda\in X'$ with $\lambda(x_k)=\delta_{jk}$ and
$\mu\in Y'$ with $\mu(y_j)\neq0$. The singleton sets
$\{\lambda\}$ and $\{\mu\}$ are equicontinuous, and
\[
 \varepsilon_{\{\lambda\},\{\mu\}}(u)=|\mu(y_j)|>0.
\]
Thus the injective seminorms separate points; the inequality already
proved also gives separation for $\pi$. The identity from
$X\otimes_\pi Y$ to $X\otimes_\varepsilon Y$ is continuous by that
inequality. Its composition with the inclusion into the injective
completion has complete Hausdorff codomain, so it extends uniquely
to the projective completion. This proves both continuity assertions.
\end{proof}

The third topology we need represents separately continuous bilinear maps. For a subset $U\subseteq X\otimes Y$ and $x\in X$, $y\in Y$, write
\[
U_x:=\{v\in Y\mid x\otimes v\in U\},
\qquad
{}_{y}U:=\{u\in X\mid u\otimes y\in U\}.
\]

\begin{definition}[Inductive tensor topology]
\label{def:topologia-tensorial-inductiva}
\index{inductive tensor product}
Let $X$ and $Y$ be locally convex spaces. The \textit{inductive tensor topology} $\iota$ on $X\otimes Y$ is the locally convex topology whose neighborhood base at the origin consists of absolutely convex absorbing subsets $U\subseteq X\otimes Y$ such that $U_x$ is a neighborhood of the origin in $Y$ for every $x\in X$ and ${}_{y}U$ is a neighborhood of the origin in $X$ for every $y\in Y$. Denote the corresponding spaces by $X\otimes_{\iota}Y$ and $X\widehat\otimes_{\iota}Y$.
\end{definition}

\begin{proposition}[Universal property of the inductive product]
\label{prop:propiedad-universal-producto-inductivo}
The topology $\iota$ is the finest locally convex topology for which the canonical map $X\times Y\to X\otimes Y$ is separately continuous. If $Z$ is locally convex and $b\colon X\times Y\longrightarrow Z$ is bilinear, then $b$ is separately continuous if and only if its linearization $\widetilde b\colon X\otimes_{\iota}Y\longrightarrow Z$ is continuous. Consequently,
\[
\varepsilon\subseteq\pi\subseteq\iota.
\]
\end{proposition}

\begin{proof}
The family described in the definition is stable under finite intersections and multiplication by nonzero scalars; moreover, for each of its elements $U$, the set $\displaystyle\frac{1}{2}U$ belongs to the same family and satisfies $\displaystyle\frac{1}{2}U+\displaystyle\frac{1}{2}U\subseteq U$. By Proposition~\ref{prop: base-discos-topologia-lc}, it determines a locally convex topology.

For $x\in X$, the preimage of $U$ under the linear map $y\mapsto x\otimes y$ is precisely $U_x$; similarly, the preimage under $x\mapsto x\otimes y$ is ${}_yU$. Thus the canonical map is separately continuous. If $\tau$ is another locally convex topology with this property and $V$ is an absolutely convex neighborhood of the origin for $\tau$, then $V_x$ and ${}_yV$ are neighborhoods for every $x$ and every $y$; hence $V$ is a neighborhood for $\iota$. This proves maximality.

Now let $b$ be separately continuous, and let $V$ be an absolutely convex neighborhood of the origin in $Z$. Then $U:=\widetilde b^{-1}(V)$ is absolutely convex and absorbing. For each $x$, $U_x$ is the preimage of $V$ under the continuous map $y\mapsto b(x,y)$; likewise, ${}_yU$ is a neighborhood. Thus $U$ is a neighborhood in $X\otimes_{\iota}Y$, and $\widetilde b$ is continuous. The converse follows by composition with the canonical map. Finally, every continuous bilinear map is separately continuous, so maximality of the two topologies gives $\pi\subseteq\iota$; the other inclusion was proved in Proposition~\ref{prop:epsilon-menos-fina-pi}.
\end{proof}

\begin{remark}[A caution for $LF$ spaces]
\label{rem:precaucion-pi-iota-LF}
Separate continuity does not generally imply joint continuity on the product of two $LF$ spaces. For this reason, the global product of two test function spaces will be described using $\widehat\otimes_{\iota}$ rather than $\widehat\otimes_{\pi}$. The projective topology will be appropriate once supports are fixed, because the factors are then Fréchet spaces and the following proposition applies.
\end{remark}

\begin{proposition}[Separate continuity on Fréchet spaces]
\label{prop:bilineal-separadamente-continua-Frechet}
\index{bilinear map!separate continuity}
Let $X$ and $Y$ be Fréchet spaces, let $Z$ be a locally convex space, and let $b\colon X\times Y\longrightarrow Z$ be a separately continuous bilinear map. Then $b$ is continuous.

More precisely, if $(p_j)_{j\in\mathbb N}$ and $(q_j)_{j\in\mathbb N}$ are increasing families of seminorms defining the topologies of $X$ and $Y$, respectively, then, for every continuous seminorm $r$ on $Z$, there exist $C>0$ and indices $k,\ell\in\mathbb N$ such that
\[
r(b(x,y))
\leq
C\,p_k(x)q_\ell(y)
\qquad
\text{for every }x\in X\text{ and every }y\in Y.
\]
\end{proposition}

\begin{proof}
Fix $r$. For $a,j\in\mathbb N$, define
\[
S_{a,j}
:=
\{x\in X\mid r(b(x,y))\leq a q_j(y)\text{ for every }y\in Y\}.
\]
For each $x\in X$, continuity of the linear map $y\mapsto b(x,y)$ provides $a$ and $j$ such that $x\in S_{a,j}$; thus $X=\displaystyle\bigcup_{a,j\in\mathbb N}S_{a,j}$. Each $S_{a,j}$ is closed: if $x_\nu\to x$ and $x_\nu\in S_{a,j}$, continuity of $x\mapsto b(x,y)$ for each $y$ allows passage to the limit in the defining inequality.

Since $X$ is a complete metric space, Baire's theorem implies that some $S_{a,j}$ has nonempty interior. The set $S_{a,j}$ is absolutely convex. If $x_0+U\subseteq S_{a,j}$ for some neighborhood $U$ of the origin, then $U-U\subseteq S_{a,j}-S_{a,j}\subseteq2S_{a,j}$; consequently, $S_{a,j}$ contains a neighborhood of the origin. Thus there exist $k\in\mathbb N$ and $\delta>0$ such that $p_k(x)<\delta$ implies $x\in S_{a,j}$. By homogeneity,
\[
r(b(x,y))
\leq
\frac{2a}{\delta}p_k(x)q_j(y).
\]
When $p_k(x)=0$, we have $p_k(tx)<\delta$ for every $t>0$; thus $tx\in S_{a,j}$ and $t\,r(b(x,y))\leq a q_j(y)$. Letting $t\to\infty$ gives $r(b(x,y))=0$. This proves the estimate with $C=\displaystyle\frac{2a}{\delta}$ and, by Proposition~\ref{prop:propiedad-universal-producto-proyectivo}, joint continuity.
\end{proof}

\subsection{Nuclear spaces}

Nuclearity allows us to compare the topologies $\pi$ and $\varepsilon$. We begin with nuclear operators between Banach spaces.

\begin{definition}[Nuclear operator]
\label{def:operador-nuclear-banach}
\index{nuclear operator}
Let $G$ and $H$ be Banach spaces. An operator $T\in\mathcal L(G,H)$ is \textit{nuclear} if there exist sequences $(\lambda_j)_{j\in\mathbb N}\subseteq\mathbb K$, $(g_j')_{j\in\mathbb N}\subseteq G'$, and $(h_j)_{j\in\mathbb N}\subseteq H$ such that
\[
\displaystyle\sum_{j=1}^{\infty}|\lambda_j|<\infty,
\qquad
\|g_j'\|_{G'}\leq1,
\qquad
\|h_j\|_H\leq1,
\]
and
\[
Tg
=
\displaystyle\sum_{j=1}^{\infty}\lambda_jg_j'(g)h_j
\qquad
\text{for every }g\in G.
\]
The series converges absolutely and uniformly for $g$ in the unit ball of $G$.
\end{definition}

\begin{lemma}[Ideal property]
\label{lem:propiedad-ideal-operadores-nucleares}
If $T\colon G\longrightarrow H$ is nuclear and $R\colon G_0\longrightarrow G$ and $S\colon H\longrightarrow H_0$ are continuous operators between Banach spaces, then $S\circ T\circ R\colon G_0\longrightarrow H_0$ is nuclear.
\end{lemma}

\begin{proof}
If $Tg=\displaystyle\sum_{j=1}^{\infty}\lambda_jg_j'(g)h_j$, then
\[
(S\circ T\circ R)g_0
=
\displaystyle\sum_{j=1}^{\infty}\lambda_j(g_j'\circ R)(g_0)S(h_j).
\]
If $R$ and $S$ are nonzero, write $g_j'\circ R=\|R\|\widetilde g_j'$ and $S(h_j)=\|S\|\widetilde h_j$, with norms at most one, and absorb $\|R\|\|S\|$ into the coefficients. The case in which either operator is zero is immediate.
\end{proof}

Let $X$ be a locally convex Hausdorff space, and let $p$ be a continuous seminorm. Denote by $X_p$ the completion of the normed space $X/\ker p$, whose norm is $\|[x]\|_p:=p(x)$, and by $Q_p\colon X\longrightarrow X_p$ the canonical map. If $q\geq p$, the identity on $X$ induces an operator of norm at most one
\[
J_{q,p}\colon X_q\longrightarrow X_p,
\qquad
J_{q,p}Q_qx:=Q_px.
\]

\begin{definition}[Nuclear space]
\label{def:espacio-localmente-convexo-nuclear}
\index{nuclear space}
A locally convex Hausdorff space $X$ is \textit{nuclear} if, for every continuous seminorm $p$ on $X$, there exists a continuous seminorm $q\geq p$ such that the operator $J_{q,p}\colon X_q\longrightarrow X_p$ is nuclear.
\end{definition}

This definition is independent of the family of seminorms used to describe the topology: if a cofinal family of seminorms satisfies the condition, every continuous seminorm is dominated, up to a positive constant, by a finite maximum of members of that family, and the ideal property allows the condition to be transferred to it.

\begin{theorem}[Nuclearity and equality of the tensor product topologies]
\label{teo:nuclearidad-pi-igual-epsilon}
Let $X$ be a nuclear space and let $Y$ be a Hausdorff locally convex space. Then the projective and injective topologies coincide on $X\otimes Y$. In particular, the algebraic identity induces a canonical topological isomorphism
\[
X\widehat\otimes_{\pi}Y
\cong
X\widehat\otimes_{\varepsilon}Y.
\]
\end{theorem}

\begin{proof}
By Proposition~\ref{prop:epsilon-menos-fina-pi}, it suffices to prove
that every projective seminorm is continuous for the injective topology.
Let $p$ and $r$ be continuous seminorms on $X$ and $Y$, respectively.
Choose $q\geq p$ so that $J_{q,p}\colon X_q\to X_p$ is nuclear,
and write
\[
 J_{q,p}z=\sum_{j=1}^\infty\lambda_j\alpha_j(z)z_j,
 \qquad
 \sum_{j=1}^\infty|\lambda_j|<\infty,
 \quad\|\alpha_j\|\leq1,\quad\|z_j\|_{X_p}\leq1.
\]
The functionals $a_j:=\alpha_jQ_q$ form an equicontinuous set
$A\subseteq X'$, since $|a_j(x)|\leq q(x)$. Let $B\subseteq Y'$
be the image of the unit ball of $Y_r'$ under composition with $Q_r$.
It is equicontinuous as well, and Hahn--Banach gives, for $v\in Y_r$,
$\|v\|_{Y_r}=\sup_{\|\eta\|\leq1}|\eta(v)|$.

Fix $u=\sum_{i=1}^sx_i\otimes y_i$ and set $\rho:=p\otimes_\pi r$.
Hahn--Banach gives a linear functional $\ell$ on $X\otimes Y$ with
\[
 |\ell(w)|\leq\rho(w)\quad(w\in X\otimes Y),
 \qquad\ell(u)=\rho(u).
\]
If $\rho(u)=0$, take $\ell=0$; otherwise extend the functional
$cu\mapsto c\rho(u)$ defined on $\mathbb Ku$.
The bilinear form $\beta(x,y):=\ell(x\otimes y)$ satisfies
$|\beta(x,y)|\leq p(x)r(y)$. It vanishes when either argument
belongs to the kernel of its respective seminorm, and hence descends
to the normed quotients. It extends to a bilinear form
$\widehat\beta\colon X_p\times Y_r\to\mathbb K$ with
\[
 |\widehat\beta(\xi,v)|\leq\|\xi\|_{X_p}\|v\|_{Y_r}.
\]
To justify the extension, for Cauchy sequences $\xi_n$ and $v_n$
in the quotients use
\[
 |\beta(\xi_n,v_n)-\beta(\xi_m,v_m)|
 \leq\|\xi_n-\xi_m\|\|v_n\|
       +\|\xi_m\|\|v_n-v_m\|.
\]
This proves existence and independence of the limit and preserves
the bilinear bound on the completions.

Define $v_j:=\sum_{i=1}^sa_j(x_i)Q_ry_i$. The dual characterization
of the norm gives
\[
 \|v_j\|_{Y_r}
 =\sup_{b\in B}\left|\sum_{i=1}^sa_j(x_i)b(y_i)\right|
 \leq\varepsilon_{A,B}(u).
\]
The nuclear representation and continuity of $\widehat\beta$ imply
\[
 \begin{aligned}
 \rho(u)
 &=\left|\sum_{i=1}^s\widehat\beta(Q_px_i,Q_ry_i)\right|\\
 &=\left|\sum_{j=1}^\infty\lambda_j
                      \widehat\beta(z_j,v_j)\right|\\
 &\leq\left(\sum_{j=1}^\infty|\lambda_j|\right)
                    \varepsilon_{A,B}(u).
 \end{aligned}
\]
The series is absolutely convergent by the last bound, and the sum
over $i$ is finite. This estimate applies directly to the seminorm
of the original tensor. Every projective seminorm is therefore
continuous for $\varepsilon$; the topologies coincide, and their
Hausdorff completions are canonically isomorphic.
\end{proof}

We shall need two stability properties. We state them only in the generality used below.

\begin{proposition}[Subspaces and finite products]
\label{prop:estabilidad-nuclear-subespacios-productos}
Every subspace of a nuclear space, equipped with the subspace topology, is nuclear. Every finite product of nuclear spaces is nuclear. Consequently, nuclearity is preserved under topological isomorphisms and on passing to complemented or closed subspaces.
\end{proposition}

Inheritance by subspaces is the theorem in
\cite[(50.3)]{Treves1967}.

\begin{proof}
The assertion about subspaces is the hereditary theorem cited above. It does not follow solely from the ideal property, since the inclusion of the subspace into the codomain need not admit a continuous projection.

Now let $X=X_1\times\cdots\times X_m$. Every continuous seminorm $p$ on
$X$ is dominated by a seminorm of the form $p_1+\cdots+p_m$, where
$p_i$ also denotes its pullback by the projection $X\to X_i$.
For each $i$ we choose $q_i\geq p_i$ so that $J_{q_i,p_i}$ is
nuclear and put $q:=q_1+\cdots+q_m$. The identifications
\[
X_q\cong\prod_{i=1}^m(X_i)_{q_i},
\qquad
X_{p_1+\cdots+p_m}\cong\prod_{i=1}^m(X_i)_{p_i}
\]
show that $J_{q,p_1+\cdots+p_m}$ is the finite diagonal sum of the
operators $J_{q_i,p_i}$. It is therefore nuclear: simply concatenate their
nuclear representations. The ideal property allows us to pass from
$p_1+\cdots+p_m$ to $p$. This proves the assertion about finite products.
Stability under topological isomorphisms follows by composing the
canonical maps with the isomorphism and its inverse; complemented
and closed subspaces are particular cases of the hereditary theorem.
\end{proof}

The basic model of a nuclear Fréchet space is the space of rapidly decreasing sequences. For $d\in\mathbb N$, let
\[
s(\mathbb Z^d)
:=
\left\{
a=(a_\nu)_{\nu\in\mathbb Z^d}
\ \middle|\
p_k(a):=\sup_{\nu\in\mathbb Z^d}(1+\|\nu\|)^k|a_\nu|<\infty
\text{ for every }k\in\mathbb N_0
\right\}.
\]
We equip it with the Fréchet topology defined by the increasing family $(p_k)_{k\in\mathbb N_0}$.

\begin{theorem}[Nuclearity of rapidly decreasing sequences]
\label{teo:sucesiones-rapidas-nuclear}
The space $s(\mathbb Z^d)$ is nuclear.
\end{theorem}

\begin{proof}
The countable increasing family $(p_k)$ separates points and defines
a metrizable locally convex topology. To check completeness, let
$(a^{(j)})$ be Cauchy for every $p_k$. Each coordinate has a limit
$a_\nu$. If $p_k(a^{(j)}-a^{(\ell)})\leq\varepsilon$ for
$j,\ell\geq j_0$, passing to the limit in each coordinate as
$\ell\to\infty$ and taking the supremum gives
$p_k(a^{(j)}-a)\leq\varepsilon$ for $j\geq j_0$.
This implies $p_k(a)<\infty$ and convergence in every seminorm.
Therefore $s(\mathbb Z^d)$ is a Fréchet space.

For $\gamma>d$, the series $\sum_{\nu\in\mathbb Z^d}(1+\|\nu\|)^{-\gamma}$
is finite. Indeed, the sets
\[
 A_j:=\{\nu\in\mathbb Z^d\mid2^j\leq1+\|\nu\|<2^{j+1}\},
 \qquad j\in\mathbb N_0,
\]
have at most $C_d2^{jd}$ elements, since they lie in a cube of
side length $2^{j+2}$. Thus
\[
 \sum_{\nu\in\mathbb Z^d}(1+\|\nu\|)^{-\gamma}
 \leq C_d\sum_{j=0}^\infty2^{j(d-\gamma)}<\infty.
\]

Fix $k\in\mathbb N_0$ and set $q:=p_{k+d+1}$.
If $e_\nu$ is the coordinate sequence, the functional
$\alpha_\nu(a):=(1+\|\nu\|)^{k+d+1}a_\nu$ extends to
$s(\mathbb Z^d)_q$ with norm at most one, while
$z_\nu:=(1+\|\nu\|)^{-k}Q_{p_k}e_\nu$ has norm one in
$s(\mathbb Z^d)_{p_k}$. For
$\lambda_\nu:=(1+\|\nu\|)^{-d-1}$, the preceding estimate gives
$\sum_\nu\lambda_\nu<\infty$. Therefore
\[
 Tz:=\sum_{\nu\in\mathbb Z^d}\lambda_\nu\alpha_\nu(z)z_\nu
\]
defines a nuclear operator on the completions: the series converges
absolutely in norm, uniformly for $\|z\|_q\leq1$.
For $a\in s(\mathbb Z^d)$, each summand is $a_\nu Q_{p_k}e_\nu$.
The truncations $a^{(R)}$ to $\|\nu\|\leq R$ satisfy
\[
 p_k(a-a^{(R)})\leq(1+R)^{-1}p_{k+1}(a)\longrightarrow0.
\]
Consequently, $TQ_qa=Q_{p_k}a$. Since $Q_q(s(\mathbb Z^d))$ is
dense in the completion, $T=J_{q,p_k}$.

Finally, every continuous seminorm $p$ is bounded above by $Cp_k$
for some $k$ and some $C\geq1$. Set $q':=Cq$.
The identity induces a bounded operator
$L\colon s(\mathbb Z^d)_{p_k}\to s(\mathbb Z^d)_p$ with
$\|L\|\leq C$. The canonical operator associated with $q'\geq p$
factors as
\[
 J_{q',p}=LJ_{q,p_k}J_{q',q}.
\]
The ideal property, Lemma~\ref{lem:propiedad-ideal-operadores-nucleares},
makes it nuclear. This verifies the definition for every continuous seminorm.
\end{proof}

\begin{theorem}[Smooth Fourier series]
\label{teo:fourier-Cinfty-toro-s}
Let $\mathbb T^d:=\bigl(\mathbb R/(2\pi\mathbb Z)\bigr)^d$. Over complex spaces, the Fourier coefficient map
\[
\mathcal F_{\mathbb T}\colon C^\infty(\mathbb T^d)\longrightarrow s(\mathbb Z^d),
\qquad
\mathcal F_{\mathbb T}(f)_\nu
:=
\frac{1}{(2\pi)^d}
\int_{\mathbb T^d}f(x)e^{-i\langle\nu,x\rangle}\,dx,
\]
is a topological isomorphism. Its inverse is given by
\[
a\longmapsto
\displaystyle\sum_{\nu\in\mathbb Z^d}a_\nu e^{i\langle\nu,x\rangle},
\]
where the series and all its differentiated series converge uniformly. In particular, $C^\infty(\mathbb T^d)$ is nuclear.
\end{theorem}

\begin{proof}
Write $q_m(f):=\max_{|\alpha|\leq m}\|D^\alpha f\|_\infty$,
where $D^\alpha$ denotes the ordinary derivative. Periodicity removes
the boundary terms in integration by parts and gives
\[
 (i\nu)^\alpha\widehat f(\nu)
 =\frac{1}{(2\pi)^d}\int_{\mathbb T^d}
       D^\alpha f(x)e^{-i\langle\nu,x\rangle}\,dx.
\]
With $\Delta=\sum_{j=1}^d\partial_j^2$, applying this identity to
$(1-\Delta)^rf$ gives
\[
 (1+\|\nu\|^2)^r|\widehat f(\nu)|\leq C_{d,r}q_{2r}(f).
\]
For each $k$, choose $2r\geq k$. Since
$(1+\|\nu\|)^k\leq2^{\displaystyle\frac{k}{2}}(1+\|\nu\|^2)^{\displaystyle\frac{k}{2}}$,
this proves $p_k(\widehat f)\leq2^{\displaystyle\frac{k}{2}}C_{d,r}q_{2r}(f)$.
Thus $\mathcal F_{\mathbb T}$ takes values in $s(\mathbb Z^d)$ and is continuous.

For $a\in s(\mathbb Z^d)$ and a multi-index $\alpha$, take
$N>|\alpha|+d$. The lattice summation estimate proved in
Theorem~\ref{teo:sucesiones-rapidas-nuclear} gives
\[
 \sum_{\nu\in\mathbb Z^d}|(i\nu)^\alpha a_\nu|
 \leq p_N(a)\sum_{\nu\in\mathbb Z^d}
                  (1+\|\nu\|)^{|\alpha|-N}<\infty.
\]
The Weierstrass test proves uniform convergence of the function
series and of every differentiated series. Its sum $g_a$ is smooth,
and $D^\alpha g_a=\sum_\nu(i\nu)^\alpha a_\nu e^{i\langle\nu,x\rangle}$:
this identity follows by passing to the limit in the integral
formula for derivatives of the finite sums, successively in each
coordinate. The same bound proves continuity of $a\mapsto g_a$.
Termwise integration, justified by uniform convergence, and
orthogonality of the exponentials give $\widehat{g_a}(\nu)=a_\nu$.

We prove uniqueness of Fourier coefficients for continuous functions.
For $N\in\mathbb N_0$ and $\theta\in\mathbb T$, let
\[
 F_N(\theta):=\frac{1}{N+1}\left|\sum_{j=0}^Ne^{ij\theta}\right|^2
 =\sum_{|\ell|\leq N}
       \left(1-\frac{|\ell|}{N+1}\right)e^{i\ell\theta}.
\]
This kernel is nonnegative and has integral one with respect to
$\displaystyle\frac{d\theta}{2\pi}$. For $0<\delta<\pi$, the geometric sum formula
implies, using the representative $\theta\in[-\pi,\pi]$,
\[
 \sup_{\delta\leq|\theta|\leq\pi}F_N(\theta)
 \leq\frac{1}{(N+1)\sin^2(\frac{\delta}{2})}.
\]
On $\mathbb T^d$, set $K_N(y):=\prod_{j=1}^dF_N(y_j)$ and
$d\mu(y):=(2\pi)^{-d}dy$. Its integral is one, and the union
of the $d$ sets $\{|y_j|\geq\delta\}$ gives
\[
 \int_{\max_j|y_j|\geq\delta}K_N(y)\,d\mu(y)
 \leq\frac{d}{(N+1)\sin^2(\frac{\delta}{2})}.
\]
Define $(g*K_N)(x):=\int_{\mathbb T^d}g(x-y)K_N(y)\,d\mu(y)$.
If $g$ is continuous, uniform continuity permits choosing $\delta$
so that $|g(x-y)-g(x)|<\varepsilon$ whenever
$\max_j|y_j|<\delta$, uniformly in $x$. Then
\[
 \|g*K_N-g\|_\infty
 \leq\varepsilon+
       \frac{2d\|g\|_\infty}{(N+1)\sin^2(\frac{\delta}{2})}.
\]
Hence $g*K_N\to g$ uniformly. Since $K_N$ is a trigonometric
polynomial, its finite expansion gives
\[
 (g*K_N)(x)=
 \sum_{\substack{\nu\in\mathbb Z^d\\|\nu_j|\leq N\ (1\leq j\leq d)}}
 \left[\prod_{j=1}^d\left(1-\frac{|\nu_j|}{N+1}\right)\right]
 \widehat g(\nu)e^{i\langle\nu,x\rangle}.
\]
If all coefficients of $g$ vanish, every convolution vanishes,
and the uniform limit implies $g=0$. Applied to
$f-g_{\mathcal F_{\mathbb T}(f)}$, this proves the other inverse
identity. Complex nuclearity follows from
Theorem~\ref{teo:sucesiones-rapidas-nuclear} and stability under
topological isomorphisms.

For real functions, the Fourier image is the closed real subspace
\[
 s_{\mathrm{sim}}:=\{a\in s(\mathbb Z^d;\mathbb C)
                   \mid a_{-\nu}=\overline{a_\nu}\}.
\]
The symmetry follows by conjugating the integral; conversely,
grouping the terms $\nu$ and $-\nu$ makes $g_a$ real.
The complex space $s(\mathbb Z^d;\mathbb C)$ remains nuclear when
regarded as a real space. In the representations of the preceding
theorem, each complex term decomposes, for $\lambda_j\neq0$, as
\[
 \lambda_j\alpha_j(z)z_j
 =|\lambda_j|\bigl(
      \operatorname{Re}\alpha_j(z)u_j
     +\operatorname{Im}\alpha_j(z)\,iu_j\bigr),
 \qquad u_j:=\frac{\lambda_j}{|\lambda_j|}z_j.
\]
The two real functionals and the two vectors have norm at most one,
and the sum of the coefficients is $2\sum_j|\lambda_j|<\infty$.
The same family of norms $(p_k)$ defines the real topology, so
the nuclearity proof remains valid. The hereditary result in
Proposition~\ref{prop:estabilidad-nuclear-subespacios-productos}
implies that $s_{\mathrm{sim}}$, and hence $C^\infty(\mathbb T^d;\mathbb R)$,
is nuclear.
\end{proof}

\begin{proposition}[Injectivity of $\varepsilon$]
\label{prop:caracter-inyectivo-epsilon}
Let $X_0\hookrightarrow X$ and $Y_0\hookrightarrow Y$ be linear topological embeddings between Hausdorff locally convex spaces. Then the canonical algebraic map
\[
X_0\otimes_{\varepsilon}Y_0
\longrightarrow
X\otimes_{\varepsilon}Y
\]
is a topological embedding. The induced map between the completions is also a topological embedding.
\end{proposition}

\begin{proof}
Identify $X_0$ and $Y_0$ with their images. Algebraic injectivity
follows from exactness of the tensor product of vector spaces over
a field. Let $A_0\subseteq X_0'$ be equicontinuous, and choose a
continuous seminorm $p_0$ with $|\lambda_0(x)|\leq p_0(x)$ for
$\lambda_0\in A_0$. The induced topology gives a continuous
seminorm $p$ on $X$ and $C>0$ with $p_0(x)\leq Cp(x)$ on $X_0$.
Hahn--Banach extends each $\lambda_0$ to a $\lambda\in X'$
with $|\lambda(x)|\leq Cp(x)$. Choosing one extension for each
element gives an equicontinuous set $A\subseteq X'$ restricting
to $A_0$. The same argument extends every equicontinuous set
$B_0\subseteq Y_0'$ to one in $Y'$.
Thus every seminorm $\varepsilon_{A_0,B_0}$ is dominated by the
restriction of a seminorm $\varepsilon_{A,B}$ of the ambient product.
Conversely, restrictions of equicontinuous sets remain equicontinuous.
This proves that the tensor topology on the subspaces agrees with
the induced topology.

Write $E:=X_0\otimes_\varepsilon Y_0$ and
$F:=X\otimes_\varepsilon Y$, and let $T\colon E\to F$ be the embedding.
The map $T$ and its inverse on $T(E)$ are uniformly continuous,
since they are linear and continuous. Let $C$ be the closure of
$T(E)$ in the completion $\widehat F$. For $z\in C$, choose a net
$Te_i\to z$. It is Cauchy in $F$; uniform continuity of the inverse
implies that $(e_i)$ is Cauchy in $E$.
Its limit in $\widehat E$ defines $Sz$. Another choice gives the
same limit, since the difference of the images converges to zero
and the embedding induces the seminorms of $E$. Passing the same
seminorm inequalities to the limit proves that
$S\colon C\to\widehat E$ is linear and continuous.
On the other hand, $T$ extends to a continuous operator
$\widehat T\colon\widehat E\to\widehat F$, whose image lies in $C$
by density of $E$. The inverse identities hold on the dense subspaces;
continuity and separation give
\[
 S\widehat T=I_{\widehat E},\qquad \widehat TS=I_C.
\]
Thus $\widehat T$ is a topological isomorphism onto $C$ and, in
particular, a topological embedding into $\widehat F$.
\end{proof}

\begin{lemma}[Smooth extension from half-spaces and orthants]
\label{lem:extension-suave-semiespacio-soporte-fijo}
For $d\in\mathbb N$ and $r\in\{0,\ldots,d\}$ put
$G_{d,r}:=\mathbb R^{d-r}\times[0,\infty)^r$. A function belongs to
$C^\infty(G_{d,r})$ if it is smooth in the interior and all its ordinary
derivatives extend continuously to each face. There is a continuous
linear operator
\[
 E_{d,r}\colon C^\infty(G_{d,r})\longrightarrow C^\infty(\mathbb R^d)
 \quad\text{such that}\quad
 (E_{d,r}f)|_{G_{d,r}}=f.
\]
Continuity refers to the derivative seminorms on compact sets.
If $K\subseteq G_{d,r}$ is compact, there is a compact set
$\widehat K\subseteq\mathbb R^d$, independent of $f$, such that
\[
 \operatorname{supp}f\subseteq K
 \ \Longrightarrow\
 \operatorname{supp}E_{d,r}f\subseteq\widehat K,
 \qquad
 \|E_{d,r}f\|_{C^m(\mathbb R^d)}
 \leq C_{m,r}\|f\|_{C^m(G_{d,r})},\quad m\in\mathbb N_0.
\]
The operators can be chosen so that extension in different groups
of variables commutes, depends continuously on smooth parameters,
and satisfies
\[
 E_{d+e,r+s}(f\boxtimes g)
 =(E_{d,r}f)\boxtimes(E_{e,s}g),
\]
after the coordinate permutation that brings the normal variables together.
\end{lemma}

\begin{proof}
We first construct coefficients that match all jets simultaneously.
For $N\in\mathbb N_0$ and $\ell\in\{0,\ldots,N\}$ define
\[
 a_{\ell,N}:=
 \prod_{\substack{j\in\{0,\ldots,N\}\\j\ne\ell}}
       \frac{1+2^j}{2^j-2^\ell}.
\]
Lagrange interpolation at the nodes $-2^0,\ldots,-2^N$,
evaluated at $1$, gives
\[
 \sum_{\ell=0}^N a_{\ell,N}(-2^\ell)^m=1,
 \qquad m\in\{0,\ldots,N\}.
\]
For each fixed $\ell$ there exists
$\displaystyle a_\ell:=\lim_{N\to\infty}a_{\ell,N}$. Indeed, for $j>\ell$ the
factor is $\displaystyle\frac{1+2^{-j}}{1-2^{\ell-j}}$, and the products of its
numerators and denominators converge to positive numbers. Let
\[
 P_+:=\prod_{j=0}^{\infty}(1+2^{-j})<\infty,
 \qquad P_-:=\prod_{j=1}^{\infty}(1-2^{-j})>0.
\]
These convergences follow by taking logarithms: apart from finitely
many factors, $|\log(1\pm2^{-j})|\leq2^{1-j}$.
Separating the factors $j<\ell$ and $j>\ell$ gives
\[
 |a_{\ell,N}|
 =2^{-\frac{\ell(\ell+1)}{2}}
 \frac{\displaystyle\prod_{\substack{0\leq j\leq N\\j\ne\ell}}
                  (1+2^{-j})}{\displaystyle\prod_{j=1}^{\ell}(1-2^{-j})
  \displaystyle\prod_{j=1}^{N-\ell}(1-2^{-j})}
 \leq\frac{P_+}{P_-^2}\,2^{-\frac{\ell(\ell+1)}{2}}.
\]
The same bound holds for $a_\ell$. In particular,
\[
 \sum_{\ell\in\mathbb N_0}|a_\ell|2^{\ell m}<\infty,
 \qquad
 \sum_{\ell\in\mathbb N_0}a_\ell(-2^\ell)^m=1,
 \qquad m\in\mathbb N_0.
\]
For the second equality, extend $a_{\ell,N}$ by zero when
$\ell>N$ and pass to the limit in the finite identity using the
preceding summable bound, with $m$ fixed and $N\geq m$.

Fix $R>0$ and $\chi\in C_c^\infty((-R,R))$ equal to one near
zero. For $f\in C^\infty(\mathbb R^{d-1}\times[0,\infty))$ put
\[
 (E_{d,1}f)(x',t):=
 \begin{cases}
 f(x',t),&t\geq0,\\
 \displaystyle\sum_{\ell\in\mathbb N_0}
   a_\ell\chi(-2^\ell t)f(x',-2^\ell t),&t<0.
 \end{cases}
\]
For each $t<0$ only finitely many indices contribute. More
precisely, if $\alpha\in\mathbb N_0^{d-1}$ and $q\in\mathbb N_0$,
the derivative of the corresponding summand is
\[
 a_\ell(-2^\ell)^q
 \bigl[\partial_\tau^q\bigl(\chi(\tau)D_{x'}^\alpha f(x',\tau)\bigr)
                      \bigr]_{\tau=-2^\ell t}.
\]
On every compact set in the variables $x'$, its supremum norm for
$t\leq0$ is bounded by
$C_q|a_\ell|2^{\ell q}$ multiplied by the largest of the norms of
the derivatives of $f$ of order at most $|\alpha|+q$ on the product
of that compact set with $[0,R]$. The Weierstrass test allows us to
differentiate the series and extend each differentiated series continuously to
$t=0$. Since $\chi=1$ near zero, its value there is
\[
 \left(\sum_{\ell\in\mathbb N_0}a_\ell(-2^\ell)^q\right)
 D_{x'}^\alpha\partial_t^q f(x',0)
 =D_{x'}^\alpha\partial_t^q f(x',0).
\]
Thus the limits of all derivatives agree from both sides.
To justify the gluing, apply the fundamental theorem of calculus
to each normal derivative across $t=0$: the primitive of the
glued continuous derivative agrees with the function on both
half-lines and at zero. By induction on the order, together with
the tangential derivatives, the glued function is $C^\infty$.
The same estimates prove continuity for the local seminorms and,
when the global norms are finite, the $C^m$ bound in the statement.

If $\operatorname{supp}f\subseteq K$, the part of the extension with
$t<0$ is supported in
$\operatorname{pr}_{x'}K\times[-R,0]$, whereas the part with
$t\geq0$ is supported in $K$. Their union is a fixed compact set.
The construction does not change the variables $x'$; when parameters
are included among them, all their derivatives satisfy the same estimates.

For $r>1$ we apply this formula successively to each normal
variable, treating the others as parameters, and for $r=0$ we take
the identity. The derivative estimates are iterated finitely many
times; the support remains in a fixed compact box.
If two different variables are extended, the two double series and
all their differentiated series converge absolutely on compact sets,
by the sums $\displaystyle \sum_{\ell\in\mathbb N_0}|a_\ell|2^{\ell m}$. Their order may therefore be interchanged.
On a product $f(x)g(y)$ each operator acts only on its own factor;
this proves the asserted commutation and product identity.
\end{proof}

\begin{corollary}[Nuclearity of test function spaces with fixed support]
\label{cor:D-K-euclidiano-nuclear}
Let $U\subseteq\mathbb R^d$ be open and let $K\subseteq U$ be compact. Then $\mathcal D_K(U)$ is a nuclear Fréchet space.
\end{corollary}

\begin{proof}
Extension by zero topologically identifies $\mathcal D_K(U)$ with the space of functions in $C^\infty(\mathbb R^d)$ whose support is contained in $K$. Choose an open cube containing $K$ whose closure is contained in the interior of a larger cube. After identifying opposite faces of the larger cube, periodic extension determines a topological embedding of $\mathcal D_K(U)$ into $C^\infty(\mathbb T^d)$. Its image is the subspace
\[
\{f\in C^\infty(\mathbb T^d)\mid f=0\text{ on }\mathbb T^d\setminus K\},
\]
which is closed because convergence in $C^\infty$ implies uniform convergence. Theorem~\ref{teo:fourier-Cinfty-toro-s} and Proposition~\ref{prop:estabilidad-nuclear-subespacios-productos} complete the proof.
\end{proof}

\begin{theorem}[Tensor products of smooth functions on tori]
\label{teo:tensorizacion-Cinfty-toros}
The map $(f,g)\mapsto((x,y)\mapsto f(x)g(y))$ induces a topological isomorphism
\[
C^\infty(\mathbb T^m)
\widehat\otimes_{\pi}
C^\infty(\mathbb T^n)
\cong
C^\infty(\mathbb T^{m+n}).
\]
\end{theorem}

\begin{proof}
The estimate
\[
 \|D_x^\alpha D_y^\beta(f\boxtimes g)\|_\infty
 \leq\|D^\alpha f\|_\infty\|D^\beta g\|_\infty
\]
proves continuity of the bilinear map.
Proposition~\ref{prop:propiedad-universal-producto-proyectivo}
extends it to the completed product. Under the Fourier isomorphism
of Theorem~\ref{teo:fourier-Cinfty-toro-s}, this operator is identified,
in the complex case, with
\[
 \Phi\colon s(\mathbb Z^m)\widehat\otimes_\pi s(\mathbb Z^n)
       \longrightarrow s(\mathbb Z^{m+n}),
 \qquad \Phi(e_\mu\otimes e_\nu)=e_{(\mu,\nu)}.
\]
For $a\in s(\mathbb Z^{m+n})$, define
\[
 R(a):=\sum_{(\mu,\nu)\in\mathbb Z^m\times\mathbb Z^n}
                a_{\mu,\nu}e_\mu\otimes e_\nu.
\]
If $N>k+\ell+m+n$, then
\[
 \begin{aligned}
 \sum_{\mu,\nu}|a_{\mu,\nu}|p_k(e_\mu)p_\ell(e_\nu)
 &\leq p_N(a)\sum_{\mu,\nu}
 \frac{(1+\|\mu\|)^k(1+\|\nu\|)^\ell}{(1+\|(\mu,\nu)\|)^N}\\
 &\leq p_N(a)\sum_{\mu,\nu}
          (1+\|(\mu,\nu)\|)^{k+\ell-N}<\infty.
 \end{aligned}
\]
The last inequality uses the lattice summation estimate on
$\mathbb Z^{m+n}$ proved in Theorem~\ref{teo:sucesiones-rapidas-nuclear}.
Thus the net of sums over finite subsets is Cauchy for every projective
seminorm and converges in the completion. The same bound proves
continuity of $R$.

Truncations of every sequence $a\in s(\mathbb Z^d)$ satisfy
\[
 p_k(a-a^{(R)})\leq(1+R)^{-1}p_{k+1}(a)\longrightarrow0.
\]
Hence sequences of finite support are dense. Their tensors also
span a dense subspace of the completed product. To verify this,
if $a^{(R)}\to a$ and $b^{(R)}\to b$, then
\[
 \begin{split}
 &(p_k\otimes_\pi p_\ell)
        (a\otimes b-a^{(R)}\otimes b^{(R)})\\
 &\quad\leq p_k(a-a^{(R)})p_\ell(b)
           +p_k(a^{(R)})p_\ell(b-b^{(R)})\longrightarrow0.
 \end{split}
\]
This approximates every elementary tensor and then every finite sum;
those sums are dense by definition of the completion.
On finite support sequences, $\Phi R=I$, and on tensors of those
sequences, $R\Phi=I$. Continuity and the two density statements
extend both identities to the complete spaces. This proves the
complex isomorphism.

For the real case, real trigonometric polynomials are dense by
symmetric Fourier truncation. In the inverse series of a real
function, group $(\mu,\nu)$ with $(-\mu,-\nu)$.
If $c=a_{\mu,\nu}$, the pair of terms is
\[
 \begin{split}
 2\operatorname{Re}\bigl(c\,e^{i\langle\mu,x\rangle}
                            e^{i\langle\nu,y\rangle}\bigr)
 &=2\operatorname{Re}c\,(C_\mu C_\nu-S_\mu S_\nu)\\
 &\quad-2\operatorname{Im}c\,(S_\mu C_\nu+C_\mu S_\nu),
 \end{split}
\]
where $C_\mu(x):=\cos\langle\mu,x\rangle$ and
$S_\mu(x):=\sin\langle\mu,x\rangle$, with the same notation for $y$.
Each summand is a product of real functions of the two sets of variables.
The zero frequency term is handled separately.
Derivative bounds for sines and cosines are controlled by the same
weights as those for exponentials. The preceding estimate, with
an additional fixed constant, proves convergence and continuity of
this inverse in the real tensor product. The same identities on
polynomials and density prove the real isomorphism.
\end{proof}

\begin{lemma}[Localized tensor products in Euclidean spaces]
\label{lem:tensorizacion-localizada-euclidiana}
Let $U\subseteq\mathbb R^m$ and $V\subseteq\mathbb R^n$ be
open, and let $K_0,K_1\subseteq U$ and $L_0,L_1\subseteq V$ be
compact sets such that
$K_0\subseteq\operatorname{Int}(K_1)$ and
$L_0\subseteq\operatorname{Int}(L_1)$. The canonical map
\[
T_{K_1,L_1}:
\mathcal D_{K_1}(U)\widehat\otimes_{\pi}\mathcal D_{L_1}(V)
\longrightarrow
\mathcal D_{K_1\times L_1}(U\times V)
\]
is a topological embedding. Its image contains $\mathcal D_{K_0\times L_0}(U\times V)$, and the inverse map restricted to this latter space is continuous.
\end{lemma}

\begin{proof}
We begin by reducing the compact sets to the fundamental cubes of two tori. Since $K_1$ and $L_1$ are bounded, there are affine diffeomorphisms $A\colon \mathbb R^m\longrightarrow\mathbb R^m$ and $B\colon \mathbb R^n\longrightarrow\mathbb R^n$ such that
\[
A(K_1)\subseteq(-\pi,\pi)^m,
\qquad
B(L_1)\subseteq(-\pi,\pi)^n.
\]
The change of variables $f\mapsto f\circ A^{-1}$ is a topological isomorphism from $\mathcal D_{K_1}(U)$ onto $\mathcal D_{A(K_1)}(A(U))$. Indeed, the chain rule expresses each derivative of $f\circ A^{-1}$ as a linear combination, with constant coefficients, of derivatives of $f$ of the same order. The inverse change satisfies the same estimates. The corresponding assertion holds for $B$ and for $A\times B$. We may therefore replace the open and compact sets by their images while keeping the original notation. We shall thus assume from now on that
\[
K_1\subseteq(-\pi,\pi)^m,
\qquad
L_1\subseteq(-\pi,\pi)^n.
\]

The inclusions $K_0\subseteq\operatorname{Int}(K_1)$ and $L_0\subseteq\operatorname{Int}(L_1)$ allow us to choose $\chi\in\mathcal D_{K_1}(U)$ and $\eta\in\mathcal D_{L_1}(V)$ equal to one on neighborhoods of $K_0$ and $L_0$, respectively. If a function in $\mathcal D_{K_1}(U)$ is extended by zero to $\mathbb R^m$, the extension vanishes in a neighborhood of the boundary of $[-\pi,\pi]^m$. It therefore determines a smooth function on $\mathbb T^m=\mathbb R^m/(2\pi\mathbb Z)^m$. We thus obtain a topological embedding
\[
J_m\colon \mathcal D_{K_1}(U)\longrightarrow C^\infty(\mathbb T^m).
\]
The embeddings $J_n$ and $J_{m+n}$ for the second factor and for the product are constructed in the same way. The topology induced by $C^\infty$ coincides with the usual topology of each space with fixed support, since in both cases it is defined by the suprema of derivatives over a fixed compact set.

In the complex case, if $\Phi\in\mathcal D_{K_0\times L_0}(U\times V)$ and
\[
\Phi(x,y)
=
\displaystyle\sum_{\mu\in\mathbb Z^m,\,\nu\in\mathbb Z^n}
c_{\mu,\nu}(\Phi)e^{i\langle\mu,x\rangle}e^{i\langle\nu,y\rangle}
\]
is its periodic Fourier series, define
\[
R_{K_0,L_0}^{K_1,L_1}\Phi
:=
\displaystyle\sum_{\substack{\mu\in\mathbb Z^m\\\nu\in\mathbb Z^n}}
c_{\mu,\nu}(\Phi)
\bigl(\chi e^{i\langle\mu,\cdot\rangle}\bigr)
\otimes
\bigl(\eta e^{i\langle\nu,\cdot\rangle}\bigr).
\]
For each $a,b\in\mathbb N_0$ there are constants $C_a,C_b>0$ such that
\[
\|\chi e^{i\langle\mu,\cdot\rangle}\|_{C^a}
\leq C_a(1+\|\mu\|)^a,
\qquad
\|\eta e^{i\langle\nu,\cdot\rangle}\|_{C^b}
\leq C_b(1+\|\nu\|)^b.
\]
By integration by parts, for every $N\in\mathbb N$ there is $C_N>0$ such that
\[
|c_{\mu,\nu}(\Phi)|
\leq
C_N(1+\|\mu\|+\|\nu\|)^{-N}
\max_{|\gamma|\leq2N}\|D^\gamma\Phi\|_{L^\infty(\mathbb T^{m+n})}.
\]
Choosing $N>a+b+m+n$, the series defining $R_{K_0,L_0}^{K_1,L_1}\Phi$ converges absolutely for the projective seminorm $\|\cdot\|_{C^a}\otimes_\pi\|\cdot\|_{C^b}$. The same estimate shows that $R_{K_0,L_0}^{K_1,L_1}$ is continuous. Since $\chi(x)\eta(y)=1$ on $K_0\times L_0$, Theorem~\ref{teo:tensorizacion-Cinfty-toros} gives
\[
T_{K_1,L_1}R_{K_0,L_0}^{K_1,L_1}\Phi
=
\chi(x)\eta(y)\Phi(x,y)
=
\Phi(x,y).
\]

In the real case choose real cutoffs $\chi$ and $\eta$.
The coefficients of a real function satisfy
$c_{-\mu,-\nu}=\overline{c_{\mu,\nu}}$.
Choose one representative of each nonzero pair
$\{(\mu,\nu),(-\mu,-\nu)\}$ and put
$C_\mu(x):=\cos\langle\mu,x\rangle$ and
$S_\mu(x):=\sin\langle\mu,x\rangle$, with the same notation for $y$.
If $c=c_{\mu,\nu}(\Phi)$, replace each pair of terms in the
inverse series by the real tensor
\[
 \begin{split}
 2\operatorname{Re}c\bigl[
    (\chi C_\mu)\otimes(\eta C_\nu)
   -(\chi S_\mu)\otimes(\eta S_\nu)\bigr]
 &\\
 -2\operatorname{Im}c\bigl[
    (\chi S_\mu)\otimes(\eta C_\nu)
   +(\chi C_\mu)\otimes(\eta S_\nu)\bigr]& .
 \end{split}
\]
The zero-frequency term is
$c_{0,0}(\Phi)\chi\otimes\eta$.
The $C^a$ and $C^b$ norms of the sine and cosine factors
satisfy the same weight bounds as the exponentials.
By the preceding estimate, the net of finite sums of these
tensors is Cauchy for every projective seminorm and converges
in the real completion. That estimate also proves continuity
of the resulting operator. Its image under $T_{K_1,L_1}$ is
the grouped Fourier series multiplied by $\chi(x)\eta(y)$,
and therefore equals $\Phi$. This constructs the localized
inverse over both scalar fields.

It remains to prove that $T_{K_1,L_1}$ induces the correct topology on its image. By Corollary~\ref{cor:D-K-euclidiano-nuclear} and Theorem~\ref{teo:nuclearidad-pi-igual-epsilon}, in the domain of $T_{K_1,L_1}$ we may replace $\widehat\otimes_\pi$ by $\widehat\otimes_\varepsilon$. The spaces $C^\infty(\mathbb T^m)$ and $C^\infty(\mathbb T^n)$ are also nuclear, so the same replacement is valid for them. Proposition~\ref{prop:caracter-inyectivo-epsilon}, applied to $J_m$ and $J_n$, shows that
\[
J_m\widehat\otimes J_n:
\mathcal D_{K_1}(U)\widehat\otimes_\pi\mathcal D_{L_1}(V)
\longrightarrow
C^\infty(\mathbb T^m)\widehat\otimes_\pi C^\infty(\mathbb T^n)
\]
is a topological embedding. If $T_{\mathbb T}$ denotes the isomorphism of Theorem~\ref{teo:tensorizacion-Cinfty-toros}, then
\[
J_{m+n}\circ T_{K_1,L_1}
=
T_{\mathbb T}\circ(J_m\widehat\otimes J_n).
\]
The right-hand side is a topological embedding. Since $J_{m+n}$ is as well, the equality proves that $T_{K_1,L_1}$ is injective and that its inverse on the image is continuous. Therefore $T_{K_1,L_1}$ is a topological embedding. Finally, the identity obtained above shows that $R_{K_0,L_0}^{K_1,L_1}$ agrees with its inverse on $\mathcal D_{K_0\times L_0}(U\times V)$.
\end{proof}

\section{Ordinary differential equations in Banach spaces}
\label{sec:edo-espacios-banach}
\label{ap:edo-banach}

The local theory of ordinary differential equations in Banach spaces retains
the method of successive approximations. Completeness of the model space
takes the place of completeness of Euclidean spaces, but no form of local
compactness is available. This distinction will be essential when studying
the continuation of integral curves on Banach manifolds.

All integrals in this section are Riemann integrals with values in a Banach
space. For continuous integrands on compact intervals, they coincide with the
corresponding Bochner integrals.

Let $E$ be a real Banach space, let $\Omega\subseteq\mathbb R\times E$ be an open subset, and let
$F\colon \Omega\longrightarrow E$. We say that $F$ is
\textbf{locally Lipschitz in the state variable, uniformly in time}
if, for each $(t_0,x_0)\in\Omega$, there exist $a,r,L>0$ such that
\[
[t_0-a,t_0+a]\times\overline B_E(x_0,r)\subseteq\Omega
\]
and
\[
\|F(t,x)-F(t,y)\|_E\leq L\|x-y\|_E
\]
for all $t\in[t_0-a,t_0+a]$ and
$x,y\in\overline B_E(x_0,r)$. When a parameter
$\lambda$ in a Banach space $P$ is present, the expression \emph{uniformly in
time and in the parameter} has the analogous meaning on products of the form
\[
[t_0-a,t_0+a]\times\overline B_E(x_0,r)
\times\overline B_P(\lambda_0,\rho).
\]

Every function of class $C^1$ is locally Lipschitz in the state variable.
Continuity of $D_2F$ allows its norm to be bounded on a sufficiently small
product around the point under consideration, and the mean value inequalities
in Theorem~\ref{teo:desigualdades-valor-medio-direccionales}
give the Lipschitz estimate on a convex ball.

\begin{theorem}[Picard--Lindelöf in Banach spaces]
\label{teo:picard-lindelof-banach}
Let $E$ be a Banach space, let $\Omega\subseteq\mathbb R\times E$ be open,
and let $F\colon \Omega\longrightarrow E$ be continuous and locally Lipschitz in
the state variable, uniformly in time. For each
$(t_0,x_0)\in\Omega$, there exist $h>0$ and a unique function
\[
x\in C^1([t_0-h,t_0+h],E)
\]
whose graph is contained in $\Omega$ and which satisfies
\[
x'(t)=F(t,x(t)),
\qquad
x(t_0)=x_0.
\]
More precisely, $a,r,L,M>0$ can be chosen so that
\[
Q:=[t_0-a,t_0+a]\times\overline B_E(x_0,r)\subseteq\Omega,
\qquad
\|F(t,y)\|_E\leq M
\]
on $Q$ and the Lipschitz constant with respect to $y$ is $L$. Any
$h>0$ such that
\[
h\leq a,
\qquad
hM\leq r,
\qquad
hL<1.
\]
then suffices.
The solution is the unique element of
$C([t_0-h,t_0+h],\overline B_E(x_0,r))$ satisfying the integral equation
\begin{equation}
\label{eq:integral-picard-banach}
x(t)=x_0+\int_{t_0}^{t}F(s,x(s))\,ds.
\end{equation}
\end{theorem}

\begin{proof}
By continuity of $F$ at $(t_0,x_0)$ and the local Lipschitz hypothesis,
we can decrease $a$ and $r$ to obtain the cylinder $Q$, the bound
$M$, and the constant $L$ in the statement. Let
\[
\mathcal X
:=
\left\{
u\in C([t_0-h,t_0+h],E)
\;\middle|\;
\sup_{t\in[t_0-h,t_0+h]}\|u(t)-x_0\|_E\leq r
\right\}.
\]
The space $C([t_0-h,t_0+h],E)$, equipped with the supremum norm, is a Banach space.
For a Cauchy sequence of continuous functions, completeness of $E$ gives
a pointwise limit; convergence is uniform, so the limit is continuous.
Consequently, $\mathcal X$ is a complete metric space, since it is a closed subset.

Define $\mathcal T\colon \mathcal X\longrightarrow C([t_0-h,t_0+h],E)$ by
\[
(\mathcal Tu)(t)
:=
x_0+\int_{t_0}^{t}F(s,u(s))\,ds.
\]
The integral is well-defined by
Theorem~\ref{teo:b5-graves-integral-riemann-banach}. For $u\in\mathcal X$,
we have
\[
\|(\mathcal Tu)(t)-x_0\|_E
\leq M|t-t_0|
\leq Mh
\leq r,
\]
so that $\mathcal T(\mathcal X)\subseteq\mathcal X$. Similarly, for
$u,v\in\mathcal X$,
\[
\|\mathcal Tu-\mathcal Tv\|_\infty
\leq hL\|u-v\|_\infty.
\]
Thus, $\mathcal T$ is a contraction.
Theorem~\ref{teo:principio-contraccion-banach} gives a unique fixed point
$x\in\mathcal X$, which satisfies
\eqref{eq:integral-picard-banach}. Since the integrand is continuous, the fundamental theorem of
calculus for Banach-valued functions,
Theorem~\ref{teo:b5-fundamental-calculo-riemann-banach}, shows that
$x\in C^1$ and $x'(t)=F(t,x(t))$.

Conversely, every solution taking values in
$\overline B_E(x_0,r)$ satisfies \eqref{eq:integral-picard-banach} by the same fundamental theorem of calculus
and is therefore a fixed point of
$\mathcal T$. If $y$ is any other solution with the same initial datum,
continuity allows both solutions to be restricted to an interval on which
they remain in a common contraction cylinder. The set of times at which
$x$ and $y$ agree is closed; repeating this local argument at each of
its points shows that it is also open. Connectedness of the common interval
shows that the two solutions agree wherever both are defined.
This proves the asserted uniqueness.
\end{proof}

The preceding uniqueness result is not restricted to the particular cylinder
used in the construction.

\begin{proposition}[Uniqueness on common domains]
\label{prop:unicidad-soluciones-edo-banach}
Under the hypotheses of Theorem~\ref{teo:picard-lindelof-banach}, let
$x_1\colon I_1\longrightarrow E$ and $x_2\colon I_2\longrightarrow E$ be two solutions,
and suppose that $t_*\in I_1\cap I_2$ and $x_1(t_*)=x_2(t_*)$. Then
\[
x_1(t)=x_2(t)
\]
for every $t$ in the connected component of $I_1\cap I_2$ containing $t_*$.
In particular, two solutions with the same initial datum agree on the
intersection of their intervals of definition.
\end{proposition}

\begin{proof}
Let $J$ be the indicated component, and let
$A:=\{t\in J\mid x_1(t)=x_2(t)\}$. This set is closed in $J$ by continuity of the solutions,
and it is nonempty. If $s\in A$, apply the local construction in
Theorem~\ref{teo:picard-lindelof-banach} at the point
$(s,x_1(s))$. By continuity, both solutions remain in the uniqueness cylinder
when restricted to a sufficiently small interval around
$s$; uniqueness of the fixed point implies that they agree on that
interval. Hence $A$ is also open in $J$.
Connectedness of $J$ implies that $A=J$.
\end{proof}

\begin{theorem}[Maximal solution]
\label{teo:solucion-maximal-edo-banach}
Under the hypotheses of Theorem~\ref{teo:picard-lindelof-banach}, for each
$(t_0,x_0)\in\Omega$ there exists a unique maximal solution
\[
x_{t_0,x_0}\colon I_{t_0,x_0}\longrightarrow E,
\]
where $I_{t_0,x_0}\subseteq\mathbb R$ is an open interval containing
$t_0$. Every other solution with initial datum $(t_0,x_0)$ is a restriction
of this solution.
\end{theorem}

\begin{proof}
Consider all solutions defined on open intervals containing
$t_0$ and taking the value $x_0$ at $t_0$. Proposition
\ref{prop:unicidad-soluciones-edo-banach} shows that any two of them agree on the intersection of their
domains. Their union therefore defines a function on the union $I_{t_0,x_0}$
of all these intervals. This union is an open interval containing
$t_0$, and the resulting function locally agrees with one of the
solutions in the family; it is therefore a solution. By construction,
it extends all the others and is maximal. Uniqueness again follows from
Proposition~\ref{prop:unicidad-soluciones-edo-banach}.
\end{proof}

The continuation criterion must be formulated in terms of convergence rather
than compactness of closed balls, since closed balls are not compact in an
infinite-dimensional Banach space.

\begin{theorem}[Continuation criterion]
\label{teo:criterio-prolongacion-edo-banach}
Under the hypotheses of Theorem~\ref{teo:picard-lindelof-banach},
let $x\colon (a,b)\longrightarrow E$ be a maximal solution of
$x'=F(t,x)$, and suppose that $b<\infty$.
\begin{enumerate}[label=(\alph*)]
\item If there exists $x_b\in E$ such that
\[
x(t)\longrightarrow x_b
\qquad\text{as }t\to b^-
\]
and $(b,x_b)\in\Omega$, then $x$ extends to an interval strictly containing
$(a,b)$, contradicting maximality.
\item If there exist $t_1\in(a,b)$ and a constant $M>0$ such that
\[
\|F(t,x(t))\|_E\leq M
\qquad
\text{for }t\in[t_1,b),
\]
then the limit $\displaystyle x_b:=\lim_{t\to b^-}x(t)$ exists. In particular, a maximal solution with a
finite right endpoint can fail to extend only if
$(b,x_b)\notin\Omega$.
\end{enumerate}
The symmetric statements hold at a finite left endpoint.
\end{theorem}

\begin{proof}
First suppose that $x(t)\to x_b$ and $(b,x_b)\in\Omega$.
Choose $a_0,r_0,L>0$ such that
\[
 [b-a_0,b+a_0]\times\overline B_E(x_b,r_0)\subseteq\Omega
\]
and $F$ is $L$-Lipschitz in the state variable on this cylinder.
By Theorem~\ref{teo:picard-lindelof-banach}, there is a solution
$y$ on an interval $I_y$ about $b$, with $y(b)=x_b$.
Convergence of $x$ to $x_b$ and continuity of $y$ allow us to
choose $h>0$ such that
\[
 [b-h,b+h]\subseteq I_y,\qquad
 [b-h,b)\subseteq(a,b),\qquad h<a_0,\qquad hL<1,
\]
and both curves take values in $\overline B_E(x_b,r_0)$ on their
respective intervals near $b$.

Define $x(b):=x_b$. The integrand $s\mapsto F(s,x(s))$ extends
continuously to $b$, with value $F(b,x_b)$.
For each $t\in[b-h,b]$, passing to the limit in the integral
equation between $t$ and a time less than $b$ gives
\[
 x(t)=x_b-\int_t^bF(s,x(s))\,ds.
\]
The curve $y$ satisfies the same identity with $y$ in place of $x$.
The number
$D:=\max_{b-h\leq t\leq b}\|x(t)-y(t)\|_E$
is finite, and the Lipschitz estimate implies
\[
 \|x(t)-y(t)\|_E
 \leq L\int_t^b\|x(s)-y(s)\|_E\,ds
 \leq hLD.
\]
Taking the maximum gives $D\leq hLD$; since $hL<1$, we obtain $D=0$.
Consequently, $x$ and $y$ agree on $[b-h,b]$, and the curve using
$x$ before $b$ and $y$ after $b$ is a solution on a strictly
larger interval. This proves \textup{(a)}.

For \textup{(b)}, if $t_1\leq s<t<b$, the fundamental theorem
of calculus gives
\[
 \|x(t)-x(s)\|_E
 \leq\int_s^t\|F(\tau,x(\tau))\|_E\,d\tau
 \leq M|t-s|.
\]
Thus $x(t)$ is Cauchy as $t\to b^-$, and completeness of
$E$ provides a limit $x_b\in E$. If $(b,x_b)\in\Omega$,
\textup{(a)} provides a continuation, contradicting maximality.
The argument at a finite left endpoint follows by reversing
the time orientation.
\end{proof}

We now record dependence on data and parameters in the form that will be used
to construct flows.

\begin{theorem}[Dependence on initial data and parameters]
\label{teo:dependencia-parametros-edo-banach}
Let $E$ and $P$ be Banach spaces, let
$\Omega\subseteq\mathbb R\times E\times P$ be open, and let
$F\colon \Omega\longrightarrow E$ be continuous and locally Lipschitz in the state variable,
uniformly in time and in the parameter. Denote by
\[
x(t;t_0,x_0,\lambda)
\]
the maximal solution of the problem
\[
x'(t)=F(t,x(t),\lambda),
\qquad
x(t_0)=x_0.
\]
Then the natural domain
\[
\mathcal D
:=
\left\{
(t,t_0,x_0,\lambda)
\;\middle|\;
t\in I_{t_0,x_0,\lambda}
\right\}
\]
is open in $\mathbb R^2\times E\times P$, and the solution map
\[
\mathcal D\longrightarrow E,
\qquad
(t,t_0,x_0,\lambda)\longmapsto x(t;t_0,x_0,\lambda),
\]
is continuous.

If $F$ is of class $C^k$, with $k\in\mathbb N\cup\{\infty\}$, the solution map is of class
$C^k$. In particular, for $k\geq1$, its derivatives with respect to the
initial state and the parameter satisfy the variational equations
\begin{align}
D_{x_0}x(t)h
&=
h+\int_{t_0}^{t}D_2F(s,x(s),\lambda)
\bigl(D_{x_0}x(s)h\bigr)\,ds,
\label{eq:variacional-dato-inicial-edo-banach}\\
D_\lambda x(t)\mu
&=
\int_{t_0}^{t}
\left[
D_2F(s,x(s),\lambda)\bigl(D_\lambda x(s)\mu\bigr)
+D_3F(s,x(s),\lambda)\mu
\right]ds.
\label{eq:variacional-parametro-edo-banach}
\end{align}
\end{theorem}

\begin{proof}
Fix initial data $(\sigma,z_0,\lambda_0)$ and denote the corresponding
solution by $x_*$. Let $J=[c,d]$ be a compact interval contained in
its maximal interval, with $\sigma\in(c,d)$, and put $T:=d-c>0$.
We first construct a common domain for nearby data.

For each $s\in J$, the hypotheses and continuity of $F$ allow us to
choose $a_s,r_s,\rho_s,L_s,M_s>0$ such that the cylinder
\[
 Q_s:=[s-a_s,s+a_s]\times\overline B_E(x_*(s),r_s)
                  \times\overline B_P(\lambda_0,\rho_s)
\]
is contained in $\Omega$, $\|F\|\leq M_s$ on $Q_s$, and
$L_s$ is a Lipschitz constant in the state variable.
Decrease $a_s$, if necessary, so that
$\|x_*(t)-x_*(s)\|_E<\displaystyle\frac{r_s}{4}$ whenever
$t\in J$ and $|t-s|<a_s$. Choose a finite subcover
$(s_i-a_i,s_i+a_i)$ of $J$ and define
\[
 \delta:=\frac{1}{4}\min_{1\leq i\leq m}r_i,\qquad
 \rho:=\frac{1}{2}\min_{1\leq i\leq m}\rho_i,\qquad
 L:=\max_{1\leq i\leq m}L_i,\qquad
 M:=\max_{1\leq i\leq m}M_i.
\]
If $t\in J$, $\|z-x_*(t)\|_E\leq\delta$, and
$\|\lambda-\lambda_0\|_P\leq\rho$, choose an index $i$
with $|t-s_i|<a_i$. Then
$\|z-x_*(s_i)\|_E<\displaystyle\frac{r_i}{2}$, and the point belongs to $Q_{s_i}$.
In particular, this entire tube lies in $\Omega$,
$\|F(t,z,\lambda)\|_E\leq M$, and the constant $L$ applies
to two states in the same ball of radius $\delta$ about
$x_*(t)$. The construction uses compactness of the graph
of $x_*$.

For $\|\lambda-\lambda_0\|_P<\rho$, put
\[
 \omega(\lambda):=
 \max_{t\in J}
 \|F(t,x_*(t),\lambda)-F(t,x_*(t),\lambda_0)\|_E.
\]
Then $\omega(\lambda)\to0$ as $\lambda\to\lambda_0$.
Indeed, for each $s\in J$ and $\varepsilon>0$, continuity
of $(t,\lambda)\mapsto F(t,x_*(t),\lambda)$ provides a
relative interval $V_s$ about $s$ and a radius $\rho_s'>0$
on which its difference from $F(s,x_*(s),\lambda_0)$ is
less than $\displaystyle\frac{\varepsilon}{2}$. A finite subcover of the $V_s$
and the minimum of their radii make the difference defining
$\omega$ less than $\varepsilon$, uniformly in $t$.
Thus choose $0<\rho_1<\rho$ such that
\[
 \omega(\lambda)<\frac{\delta}{4Te^{LT}}
 \quad\text{if }\|\lambda-\lambda_0\|_P<\rho_1,
 \qquad
 \eta:=\frac{\delta}{4e^{LT}}.
\]

Let $\|z-z_0\|_E<\eta$ and
$\|\lambda-\lambda_0\|_P<\rho_1$, and let $u$ be the
maximal solution with $u(\sigma)=z$. As long as $u$ stays
in the tube, subtracting the integral equations and applying
Lemma~\ref{lema: gronwall} gives, in both directions from $\sigma$,
\[
 \|u(t)-x_*(t)\|_E
 \leq e^{L|t-\sigma|}
       \bigl(\|z-z_0\|_E+|t-\sigma|\omega(\lambda)\bigr)
 <\frac{\delta}{2}.
\]
Write $I_u=(a_u,b_u)$ and define
\[
 \tau_+:=\sup\{q\in[\sigma,d]\cap I_u\mid
 \|u(s)-x_*(s)\|_E<\delta
 \text{ for every }s\in[\sigma,q]\}.
\]
Continuity and the initial inequality give $\tau_+>\sigma$.
If $\tau_+<\min\{b_u,d\}$, the preceding estimate and
continuity give
$\|u(\tau_+)-x_*(\tau_+)\|_E\leq\displaystyle\frac{\delta}{2}$; the distance
remains less than $\delta$ on an additional interval,
contradicting the definition of $\tau_+$. Hence
$\tau_+=\min\{b_u,d\}$. If $b_u\leq d$, the bound $M$ on
the field and Theorem~\ref{teo:criterio-prolongacion-edo-banach}
give a limit $u_b$ as $t\to b_u^-$. The estimate implies
$\|u_b-x_*(b_u)\|_E\leq\displaystyle\frac{\delta}{2}$, so
$(b_u,u_b,\lambda)\in\Omega$. The same theorem continues
$u$, contradicting maximality. Thus $b_u>d$.
Applying this argument backwards gives $a_u<c$. Consequently,
$u$ is defined throughout $J$, and
\[
 \|u-x_*\|_{C(J,E)}
 \leq e^{LT}\bigl(\|z-z_0\|_E+T\omega(\lambda)\bigr).
\]
This proves the common domain and continuous dependence
in the uniform norm.

To vary the initial time, use the field
\[
 \widehat F(s,z,\tau,\lambda):=F(\tau+s,z,\lambda)
\]
on the open set
$\widehat\Omega:=\{(s,z,\tau,\lambda)\mid
(\tau+s,z,\lambda)\in\Omega\}$. This field satisfies the
same hypotheses, with $(\tau,\lambda)$ as parameter and
fixed initial time $s=0$. Apply the preceding construction
to an interval $[\alpha,\beta]$ containing both $0$ and
$t-\sigma$ in its interior. For some positive radii
$\eta,\rho_t,\rho_p$, all points of the open set
\[
 \begin{split}
 \mathcal W:=\{(t',\tau,z,\lambda)\mid&
 |\tau-\sigma|<\rho_t,\quad \|z-z_0\|_E<\eta,\\
 &\|\lambda-\lambda_0\|_P<\rho_p,\quad
 \alpha<t'-\tau<\beta\}
 \end{split}
\]
belong to $\mathcal D$. This proves that $\mathcal D$ is
open. The uniform estimate and continuity of the reference
curve prove joint continuity upon evaluation at $s=t'-\tau$.

Now suppose $F\in C^k$. Work with $\widehat F$ so that
the initial time is included among the parameters from the
outset, and write $\theta=(\tau,\lambda)$, using the sum norm
on $\mathbb R\times P$. On an interval
$J_0=[s_0-h,s_0+h]$ about a point of a reference solution
$v_*$, choose an open tube as above on which
$\|D_z\widehat F\|\leq L_0$. Take $hL_0<\displaystyle\frac{1}{2}$.
On the open neighborhood of curves and parameters whose
graphs stay in this tube, define
\[
 \mathcal N(v,\theta)(s):=\widehat F(s,v(s),\theta).
\]
This map into $C(J_0,E)$ is $C^k$, with
\[
 \begin{split}
 D^j\mathcal N(v,\theta)
 &[(w_1,\mu_1),\ldots,(w_j,\mu_j)](s)\\
 &=D_{(z,\theta)}^j\widehat F(s,v(s),\theta)
       [(w_1(s),\mu_1),\ldots,(w_j(s),\mu_j)] .
 \end{split}
\]
We specify the uniformity needed for this assertion.
For each fixed pair $(v,\theta)$, each $j\leq k$, and each
$\varepsilon>0$, continuity of the derivative on its compact
graph provides a finite subcover and a minimum $\kappa>0$
of radii such that
\[
 \sup_{s\in J_0}
 \|D_{(z,\theta)}^j\widehat F(s,v(s)+w(s),\theta+\mu)
       -D_{(z,\theta)}^j\widehat F(s,v(s),\theta)\|
 <\varepsilon
\]
whenever $\|w\|_\infty+\|\mu\|<\kappa$. The construction
is the same as the one used for $\omega$, applied to this
derivative. It also provides a uniform bound on the derivative
norms on a sufficiently small tube.
For the first derivative the remainder is bounded by
\[
 \begin{split}
 &\|\mathcal N(v+w,\theta+\mu)-\mathcal N(v,\theta)
             -D\mathcal N(v,\theta)[w,\mu]\|_\infty\\
 &\quad\leq(\|w\|_\infty+\|\mu\|)
 \sup_{\substack{s\in J_0\\0\leq q\leq1}}
 \|D_{(z,\theta)}\widehat F(s,v(s)+qw(s),\theta+q\mu)
          -D_{(z,\theta)}\widehat F(s,v(s),\theta)\|.
 \end{split}
\]
The supremum tends to zero. Applying the same argument to
successive derivatives proves the formula and continuity
through order $k$. For $k=\infty$ this argument is carried
out for each finite order.

Now consider
\[
 \mathscr F(v,a,\theta)(s)
 :=v(s)-a-\int_{s_0}^s\widehat F(q,v(q),\theta)\,dq .
\]
Continuity and linearity of integration show that
$\mathscr F$ is $C^k$. Its derivative with respect to $v$
is $I-\mathcal K$, where
\[
 (\mathcal Kw)(s)
 :=\int_{s_0}^sD_z\widehat F(q,v(q),\theta)w(q)\,dq,
 \qquad \|\mathcal K\|\leq hL_0<\frac{1}{2} .
\]
The Neumann series in Theorem~\ref{teo:serie-de-neumann}
converges in norm and gives
\[
 (I-\mathcal K)^{-1}=\sum_{j=0}^\infty\mathcal K^j,
 \qquad \|(I-\mathcal K)^{-1}\|\leq2.
\]
Thus
Theorem~\ref{teo:funcion-implicita-banach} gives a solution
that is $C^k$ in $(a,\theta)$ with values in $C(J_0,E)$.
Picard uniqueness identifies this solution with the
maximal solution already constructed.

Differentiating the integral equation with respect to
the parameters gives continuous curves satisfying the
differentiated integral equations. The fundamental theorem
of calculus allows differentiation of these curves with
respect to $s$. The identity
$\partial_sv=\widehat F(s,v,\theta)$ and the chain rule give,
by induction, all mixed derivatives of total order at most
$k$, with joint continuity.
To cover an arbitrary compact interval, take a finite
partition with step less than $\displaystyle\frac{1}{2L_0}$, increasing
$L_0>0$ so that it bounds the state derivatives on a tube
about the whole graph. The solution is obtained by composing
the local maps at the partition points and evaluating the
last piece at the variable time. These compositions are
$C^k$. Finally, the change $s=t-\tau$ preserves this
regularity and proves the assertion on $\mathcal D$.

With the initial time fixed, differentiating the integral
equation with respect to $z_0$ and $\lambda$ gives
\eqref{eq:variacional-dato-inicial-edo-banach} and
\eqref{eq:variacional-parametro-edo-banach}.
\end{proof}

\begin{corollary}[Global existence for bounded fields]
\label{cor:existencia-global-campo-acotado-edo-banach}
Let $E$ be a Banach space, and let $F\colon \mathbb R\times E\longrightarrow E$
be continuous and locally Lipschitz in the state variable,
uniformly in time. Suppose that there exists $M\geq0$ such that
\[
\|F(t,x)\|_E\leq M
\]
for every $(t,x)\in\mathbb R\times E$. Then every maximal solution is defined on
$\mathbb R$ and satisfies
\[
\|x(t)-x(s)\|_E\leq M|t-s|
\]
for all $s,t\in\mathbb R$. If $F$ has $C^k$ regularity,
the regularity of the solution map is that given by
Theorem~\ref{teo:dependencia-parametros-edo-banach}.
\end{corollary}

\begin{proof}
Let $(a,b)$ be the maximal interval of a solution. If $b<\infty$, the bound on
$F$ and Theorem~\ref{teo:criterio-prolongacion-edo-banach} give a limit
$x_b\in E$ as $t\to b^-$. Since the domain of the field is all of
$\mathbb R\times E$, the datum $(b,x_b)$ belongs to that domain and the solution extends,
contradicting maximality. Hence $b=+\infty$.
The same argument backward in time shows that $a=-\infty$.
The estimate follows by integrating the equation between
$s$ and $t$.
\end{proof}

\chapter{Linear algebra}\label{ap:algebra-lineal}

This appendix collects some results from linear algebra used in the study of metrics, bilinear forms, operators, and tensors. We mainly follow the references \cite{Horn,Roman}; the reader may consult them for a more extensive treatment.

\section{Bilinear forms and matrices}

Bilinear forms and their associated matrices provide the starting point for the constructions considered in this appendix.
\begin{definition}\label{def:c-algebra-lineal-espacio-vectorial-real-dimension-forma-bilineal}\index{bilinear form}
 Let $V$ be a real vector space of dimension $n$. We say that $f\colon V\times V\longrightarrow\mathbb{R}$ is a bilinear form on $V$ if:
 \begin{enumerate}[label=(\alph*)]
 \item $f(x+y,z)=f(x,z)+f(y,z)$ and $f(\lambda x,y)=\lambda f(x,y)$ for any $x,y,z\in V$ and $\lambda\in\mathbb{R}$.
 \item $f(x,y+z)=f(x,y)+f(x,z)$ and $f(x,\lambda y)=\lambda f(x,y)$ for any $x,y,z\in V$ and $\lambda\in\mathbb{R}$.
 \end{enumerate}
 If $\{e_{1},\dots,e_{n}\}$ is a basis of $V$, we define the matrix associated with the bilinear form $f$ by \[A_{ij}=f(e_{i},e_{j}).\]

\end{definition}
\begin{proposition}\label{prop:c-algebra-lineal-espacio-vectorial-real-dimension-forma-bilineal}
 If $V$ is a real vector space of dimension $n$, $f\colon V\times V\longrightarrow\mathbb{R}$ is a bilinear form, and $\beta=\{e_{1},\dots,e_{n}\}$ is a basis of $V$, then \[f(x,y)=[x]_{\beta}^{T}A[y ]_{\beta},\] for any $x,y\in V$, where $A$ is the matrix associated with $f$ and $[x]_{\beta},[y]_{\beta}$ are the coordinate representations of $x$ and $y$, respectively, in the basis $\beta$.
\end{proposition}
\begin{proof}
If \(x=x^ie_i\) and \(y=y^je_j\), bilinearity gives
\[
f(x,y)=x^iy^jf(e_i,e_j)=x^iA_{ij}y^j=[x]_\beta^TA[y]_\beta.
\]
\end{proof}
\begin{definition}\label{def:c-algebra-lineal-espacio-vectorial-real-dimension-finita-forma}\index{bilinear form!symmetric}\index{bilinear form!nondegenerate}
 Let $V$ be a finite-dimensional real vector space of dimension $n$, and let $f\colon V\times V\longrightarrow\mathbb{R}$ be a bilinear form. We say that it is:
 \begin{enumerate}[label=(\alph*)]
 \item Symmetric if $f(x,y)=f(y,x)$ for any $x,y\in V$. In that case, the associated matrix $A$ is symmetric.
 \item Nondegenerate if $f(x,y)=0$ for every $y\in V$ implies $x=0$.

 \end{enumerate}
 \end{definition}
 \begin{example}\label{hessiana}
 Let $f\colon U\subseteq\mathbb{R}^{n}\longrightarrow\mathbb{R}$ be of class $C^{2}(U)$. We define the Hessian matrix of $f$ at $x\in U$ by $H(x)=\left(\displaystyle\frac{\partial ^{2}f}{\partial x_{i}\partial x_{j}}(x)\right)\in M_{n\times n}(\mathbb{R})$. Then $H(x)$ is a symmetric matrix and is associated with the bilinear form $\nabla^{2}f(x)\colon \mathbb{R}^{n}\times\mathbb{R}^{n}\longrightarrow\mathbb{R}$ (the second total derivative of $f$).
 \end{example}
 \begin{definition}\label{def:c-algebra-lineal-matriz-hermitiana-matriz-obtenida-conjugar-entradas}\index{Hermitian matrix@Hermitian matrix}
 A matrix $A\in M_{n\times n}(\mathbb{C})$ is Hermitian if $A=A^{*}$, where $A^{*}$ is the matrix obtained by conjugating the entries of the transpose of $A$. A Hermitian matrix is \emph{positive semidefinite} if $x^*Ax\geq0$ for every $x\in\mathbb C^n$, and \emph{positive definite} if $x^*Ax>0$ for every $x\neq0$. For real matrices, we use these terms for symmetric matrices satisfying the corresponding condition on $x\in\mathbb R^n$.
 \end{definition}
 \begin{proposition}\label{prop:c-algebra-lineal-matriz-simetrica-hermitiana}
 Every symmetric matrix $A\in M_{n\times n}(\mathbb{R})$ is Hermitian.
 \end{proposition}
\begin{proof}
Since the entries are real, conjugation leaves them unchanged. Therefore,
\(A^*=\overline{A^T}=A^T=A\).
\end{proof}
 \begin{theorem}\label{diagonalizable}\index{diagonalizable}
 Every $n\times n$ Hermitian matrix is diagonalizable by a unitary matrix and has real eigenvalues.
 \end{theorem}
\begin{proof}
We proceed by induction on \(n\). The characteristic polynomial of a
complex matrix has a root, so \(A\) has an eigenvector
\(v\ne0\), say \(Av=\lambda v\). Since \(A=A^*\),
\[
\lambda\langle v,v\rangle
=\langle Av,v\rangle
=\langle v,Av\rangle
=\overline\lambda\langle v,v\rangle,
\]
and \(\lambda\in\mathbb R\). After normalizing \(v\), its orthogonal
complement is invariant: if \(w\perp v\), then
\[
\langle Aw,v\rangle=\langle w,Av\rangle
=\lambda\langle w,v\rangle=0.
\]
The restriction of \(A\) to \(v^\perp\) remains Hermitian. The induction
hypothesis gives an orthonormal basis of eigenvectors in
\(v^\perp\), which together with \(v\) diagonalizes \(A\).
\end{proof}

 \begin{theorem}[Rayleigh]\label{rayleigh}\index{Rayleigh theorem@Rayleigh theorem}Let $A\in M_{n\times n}(\mathbb{C})$ be Hermitian. Then \[\lambda_{\max}=\max_{x\neq 0}\frac{x^{*}Ax}{x^{*}x} \text{ and }\lambda_{\min}=\min_{x\neq 0}\frac{x^{*}Ax}{x^{*}x},\] where $\lambda_{\displaystyle\max}$ and $\lambda_{\min}$ are the largest and smallest eigenvalues of $A$, respectively, and $x^{*}$ is the matrix obtained by transposing the vector $x$ and taking the complex conjugate of each entry.

 \end{theorem}
\begin{proof}
By the preceding theorem, there is a unitary matrix \(Q\) such that
\[
Q^*AQ=\operatorname{diag}(\lambda_1,\dots,\lambda_n).
\]
If \(y=Q^*x\), then
\[
\frac{x^*Ax}{x^*x}
=\frac{\sum_{i=1}^{n}\lambda_i|y_i|^2}{\sum_{i=1}^{n}|y_i|^2}.
\]
This quotient is a convex combination of the eigenvalues, so
it lies between \(\lambda_{\min}\) and \(\lambda_{\displaystyle\max}\). Both extrema are
attained by taking a corresponding eigenvector.
\end{proof}
 \begin{theorem}[Schur decomposition]\label{schur}\index{Schur decomposition@Schur decomposition}Let $A\in M_{n\times n}(\mathbb{C})$. Then $A=QUQ^{-1}$, where $Q$ is a unitary matrix ($Q$ is invertible and $Q^{-1}=Q^{*}$) and $U$ is an upper triangular matrix whose diagonal entries are the eigenvalues of $A$.

 \end{theorem}
\begin{proof}
The proof is by induction on \(n\). Choose a unit eigenvector
\(v_1\) of \(A\) and extend it to an orthonormal basis. In that basis, the
first column of the matrix of \(A\) is
\((\lambda_1,0,\dots,0)^T\), so the matrix has the form
\[
\begin{pmatrix}
\lambda_1& *\\
0&A_1
\end{pmatrix}.
\]
By the induction hypothesis, a unitary change of basis in
\(\operatorname{span}\{v_1\}^\perp\) triangularizes \(A_1\). Extending it
by the identity on \(\operatorname{span}\{v_1\}\) gives a
unitary matrix \(Q\) for which \(Q^*AQ\) is upper triangular. The
diagonal entries are its eigenvalues because the characteristic
polynomial of a triangular matrix is the product of
\(\lambda-U_{ii}\).
\end{proof}
 \begin{theorem}\label{valor propio y frobenius}\index{Frobenius norm@Frobenius norm!eigenvalue bound}
Let $A\in M_{n\times n}(\mathbb{C})$, and let $\lambda$ be an eigenvalue of $A$. If
\[
\|A\|_{F}:=\left(\displaystyle\sum_{i,j=1}^{n}|A_{ij}|^{2}\right)^{\frac{1}{2}}
\]
is the Frobenius norm of $A$, then $|\lambda|\leq\|A\|_{F}$.
\end{theorem}

\begin{proof}
By Theorem~\ref{schur}, there exist a unitary matrix $Q$ and an upper triangular matrix $U$ such that $A=QUQ^{-1}$. The identity $\|C\|_F^2=\operatorname{tr}(C^*C)$ and the equality $\operatorname{tr}(XY)=\operatorname{tr}(YX)$ give $\|Q^*AQ\|_F^2=\operatorname{tr}(Q^*A^*AQ)=\operatorname{tr}(A^*A)$; hence $\|A\|_{F}=\|U\|_{F}$. Since the diagonal entries of $U$ are the eigenvalues of $A$, we have
\[
\|U\|_{F}^{2}
=
\displaystyle\sum_{i=1}^{n}|U_{ii}|^{2}
+
\displaystyle\sum_{1\leq i<j\leq n}|U_{ij}|^{2}
\geq
|\lambda|^{2}.
\]
Therefore, $|\lambda|\leq\|U\|_{F}=\|A\|_{F}$.
\end{proof}

 \begin{theorem}\label{teo:c-algebra-lineal-matriz-hermitiana-positivo-semidefinida-label-unica}\index{positive root of a positive semidefinite Hermitian matrix@positive root of a positive semidefinite Hermitian matrix}
 Let $A$ be an $n\times n$ positive semidefinite Hermitian matrix, and let $k\in \mathbb{N}$, with $k\geq 2$. Then:
 \begin{enumerate}[label=(\alph*)]
 \item There exists a unique positive semidefinite Hermitian matrix $B$ such that $A=B^{k}$.
 \item There exists a polynomial $p$ with real coefficients such that $B=p(A)$. Consequently, $B$ commutes with every matrix that commutes with $A$.
 \item If $A$ is real, then $B$ is real.
 \end{enumerate}

 \end{theorem}
\begin{proof}
By the spectral theorem,
\[
A=Q\operatorname{diag}(\lambda_1,\dots,\lambda_n)Q^*,
\qquad \lambda_i\geq0,
\]
with \(Q\) unitary. Define
\[
B:=Q\operatorname{diag}
(\lambda_1^{1/k},\dots,\lambda_n^{1/k})Q^*.
\]
Then \(B\) is positive semidefinite Hermitian and \(B^k=A\).

If \(C\) is another positive semidefinite Hermitian matrix with \(C^k=A\),
then \(C\) commutes with \(A=C^k\) and preserves each eigenspace
\(E_\lambda\) of \(A\). The restriction of \(C\) to \(E_\lambda\) is
positive semidefinite Hermitian and satisfies
\((C\restriction_{E_\lambda})^k=\lambda I\). Its eigenvalues are
nonnegative and their \(k\)th powers equal \(\lambda\); hence they all equal
\(\lambda^{1/k}\). The spectral theorem implies
\(C\restriction_{E_\lambda}=\lambda^{1/k}I\), so \(C=B\).

Let \(\mu_1,\dots,\mu_s\) be the distinct eigenvalues of \(A\).
Lagrange interpolation gives a real polynomial \(p\) with
\(p(\mu_j)=\mu_j^{1/k}\). Calculation on each eigenspace shows that
\(B=p(A)\). If a matrix commutes with \(A\), it commutes with \(p(A)\); and if
\(A\) is real, so is \(p(A)\).
\end{proof}
 \begin{theorem}\label{teo:c-algebra-lineal-matriz-hermitiana-positivo-semidefinida-positivo-definida}\index{Hermitian matrix@Hermitian matrix!positive definiteness criterion}
 If $A$ is a positive semidefinite Hermitian matrix, then $A$ is positive definite if and only if $A$ is invertible.
 \end{theorem}
\begin{proof}
By the spectral theorem, \(A\) is unitarily similar to
\(\operatorname{diag}(\lambda_1,\dots,\lambda_n)\) with
\(\lambda_i\geq0\). It is positive definite if and only if all the
\(\lambda_i\) are strictly positive, and it is invertible if and only if
none of them is zero. The two conditions are equivalent.
\end{proof}
 \begin{remark}\label{raiz cuadrada} For every real symmetric positive definite matrix $A$, there exists a unique real symmetric positive definite matrix $B$ such that $A=B^{2}$. This matrix is called the square root of $A$.

 \end{remark}
 \begin{proposition}\label{raiz cuadrada simetrica}
 A symmetric positive (semi)definite matrix has a positive (semi)definite square root that is also symmetric.
 \end{proposition}
\begin{proof}
Apply Theorem
\ref{teo:c-algebra-lineal-matriz-hermitiana-positivo-semidefinida-label-unica}
with \(k=2\). Since the original matrix is real, part (c) of that theorem
shows that its positive root is real; being Hermitian and real, it is
symmetric. If the original matrix is positive definite, the root has
no zero eigenvalues and is also positive definite.
\end{proof}

\chapter{Einstein summation convention}\label{ap:convenio-einstein}

This appendix specifies the summation convention used in expressions with indices.

\section{Einstein summation convention}

In a tensor expression, an index that appears exactly twice in the same monomial, once as a subscript and once as a superscript, is called a \emph{dummy index}. Unless otherwise stated, summation over all possible values of that index is understood. Thus, if the indices range over $\{1,\dots,n\}$, we write
\[
\lambda_i x^i
:=
\displaystyle\sum_{i=1}^{n}\lambda_i x^i
\]
and
\[
g^{ij}F_iG_j
:=
\displaystyle\sum_{i,j=1}^{n}g^{ij}F_iG_j.
\]

Indices that appear only once are called \emph{free indices}. In a tensor equality, the free indices must appear on both sides in the same positions. For example,
\[
T^i{}_jX^j=Y^i
\]
represents a family of equalities, one for each value of the free index $i$, and contains a sum over the dummy index $j$.

In each monomial we will avoid using the same index more than twice. When several sums are needed, we will use distinct dummy indices. This rule avoids ambiguities and makes the free indices of the expression immediately recognizable.

Local coordinates will usually be written as $(x^1,\dots,x^n)$. In particular,
\[
\frac{\partial f}{\partial x^i}\,\mathbf{d}x^i
:=
\displaystyle\sum_{i=1}^{n}\frac{\partial f}{\partial x^i}\,\mathbf{d}x^i.
\]

\cleardoublepage
\printglossary[title={Glossary}]

\cleardoublepage
\chapter*{Image credits}
\addcontentsline{toc}{chapter}{Image credits}
\markboth{Image credits}{Image credits}

The following historical images are from Wikimedia Commons. The rights
information reproduces the entry consulted for each file; where an image
is identified as being in the public domain only in the United States,
reuse in other countries may depend on the applicable law.

\begin{itemize}
  \item \textbf{Isaac Newton.} James Thornhill's copy of the portrait by
  Godfrey Kneller (1689),
  \href{https://commons.wikimedia.org/wiki/File:GodfreyKneller-IsaacNewton-1689.jpg}{Wikimedia Commons}, public domain.
  Appears in Figures~\ref{fig:introduccion-calculo-geometria-global}
  and~\ref{fig:newton-leibniz}.

  \item \textbf{Gottfried Wilhelm von Leibniz.} Portrait by Christoph
  Bernhard Francke, ca.~1695,
  \href{https://commons.wikimedia.org/wiki/File:Gottfried_Wilhelm_von_Leibniz.jpg}{Wikimedia Commons}, public domain.
  Appears in Figures~\ref{fig:introduccion-calculo-geometria-global}
  and~\ref{fig:newton-leibniz}.

  \item \textbf{Carl Friedrich Gauss.} Portrait by Christian Albrecht
  Jensen (1840),
  \href{https://commons.wikimedia.org/wiki/File:Carl_Friedrich_Gauss_1840_by_Jensen.jpg}{Wikimedia Commons}, public domain.
  Appears in Figures~\ref{fig:introduccion-calculo-geometria-global}
  and~\ref{fig:gauss-riemann-poincare}.

  \item \textbf{Bernhard Riemann.} Photograph from 1863, unidentified
  photographer,
  \href{https://commons.wikimedia.org/wiki/File:Georg_Friedrich_Bernhard_Riemann.jpeg}{Wikimedia Commons}, public domain according to the file entry
  (\emph{PD-old}). Appears in
  Figures~\ref{fig:introduccion-calculo-geometria-global}
  and~\ref{fig:gauss-riemann-poincare}.

  \item \textbf{Henri Poincaré.} Image published as the frontispiece of
  \emph{Last Thoughts} (1913), unidentified creator,
  \href{https://commons.wikimedia.org/wiki/File:JH_Poincare.jpg}{Wikimedia Commons}, public domain in the United States
  (\emph{PD-US}). Appears in
  Figures~\ref{fig:introduccion-calculo-geometria-global}
  and~\ref{fig:gauss-riemann-poincare}.

  \item \textbf{Élie Cartan.} Photograph from 1904, unidentified photographer,
  reproduced on
  \href{https://commons.wikimedia.org/wiki/File:Elie-Cartan-1904.png}{Wikimedia Commons}, public domain in the United States
  (\emph{PD-US}). Appears in Figure~\ref{fig:elie-cartan}.

  \item \textbf{Oliver Heaviside.} Portrait by an unidentified artist,
  \href{https://commons.wikimedia.org/wiki/File:Oliver_Heaviside2.jpg}{Wikimedia Commons}, reproduction of a work in the public domain
  (\emph{PD-Art}). Appears in Figure~\ref{fig:heaviside}.
\end{itemize}

\cleardoublepage
\phantomsection
\addcontentsline{toc}{chapter}{Bibliography}
\nocite{*}
\bibliographystyle{plain}
\bibliography{bibliografia}
\index{differential operator!formal adjoint|see{formal adjoint}}
\index{Bochner, Laplacian|see{Bochner Laplacian}}
\index{differential|seealso{derivative}}
\index{Frechet, derivative@Fréchet, derivative|see{Fréchet derivative}}
\index{Gateaux, derivative@Gâteaux, derivative|see{Gâteaux derivative}}
\index{Hodge--de Rham, Laplacian|see{Hodge--de Rham Laplacian}}
\index{Laplace--Beltrami, operator|see{Laplace--Beltrami operator}}
\index{Nehari, manifold|see{manifold, Nehari}}
\index{Riemann--Lebesgue, measure|see{Riemann--Lebesgue measure}}
\index{Sobolev, space|see{Sobolev space}}
\index{derivative|seealso{differential}}
\index{method|seealso{variational principle}}
\index{differential operator|seealso{principal symbol}}
\index{manifold|seealso{tangent space}}
\index{Hilbert, manifold|see{Hilbert manifold}}
\index{Finsler, Banach-- structure|see{Banach--Finsler structure}}
\index{Palais--Smale, condition|see{Palais--Smale condition}}
\index{Eichhorn, uniformity|see{uniform structure, Eichhorn}}
\index{Amann, uniform regularity|see{Riemannian manifold, uniformly regular}}
\index{ILH, manifold|see{ILH manifold}}

\printindex

\end{document}